\documentclass[numbers=enddot,12pt,final,onecolumn,notitlepage]{scrartcl}%
\usepackage[headsepline,footsepline,manualmark]{scrlayer-scrpage}
\usepackage[all,cmtip]{xy}
\usepackage{amssymb}
\usepackage{amsmath}
\usepackage{amsthm}
\usepackage{framed}
\usepackage{comment}
\usepackage{color}
\usepackage[breaklinks=True]{hyperref}
\usepackage[sc]{mathpazo}
\usepackage[T1]{fontenc}
\usepackage{tikz}
\usepackage{needspace}
\usepackage{tabls}
\usepackage{hyperxmp}
\usepackage[type={CC}, modifier={zero}, version={1.0},]{doclicense}
\providecommand{\U}[1]{\protect\rule{.1in}{.1in}}
\usetikzlibrary{arrows}
\newcounter{exer}

\numberwithin{exer}{section}
\theoremstyle{definition}
\newtheorem{theo}{Theorem}[section]
\newenvironment{theorem}[1][]
{\begin{theo}[#1]\begin{leftbar}}
{\end{leftbar}\end{theo}}
\newtheorem{lem}[theo]{Lemma}
\newenvironment{lemma}[1][]
{\begin{lem}[#1]\begin{leftbar}}
{\end{leftbar}\end{lem}}
\newtheorem{prop}[theo]{Proposition}
\newenvironment{proposition}[1][]
{\begin{prop}[#1]\begin{leftbar}}
{\end{leftbar}\end{prop}}
\newtheorem{defi}[theo]{Definition}
\newenvironment{definition}[1][]
{\begin{defi}[#1]\begin{leftbar}}
{\end{leftbar}\end{defi}}
\newtheorem{remk}[theo]{Remark}
\newenvironment{remark}[1][]
{\begin{remk}[#1]\begin{leftbar}}
{\end{leftbar}\end{remk}}
\newtheorem{coro}[theo]{Corollary}
\newenvironment{corollary}[1][]
{\begin{coro}[#1]\begin{leftbar}}
{\end{leftbar}\end{coro}}
\newtheorem{conv}[theo]{Convention}
\newenvironment{convention}[1][]
{\begin{conv}[#1]\begin{leftbar}}
{\end{leftbar}\end{conv}}
\newtheorem{quest}[theo]{Question}

\newtheorem{warn}[theo]{Warning}

\newtheorem{conj}[theo]{Conjecture}

\newtheorem{exam}[theo]{Example}
\newenvironment{example}[1][]
{\begin{exam}[#1]\begin{leftbar}}
{\end{leftbar}\end{exam}}
\newtheorem{exmp}[exer]{Exercise}
\newenvironment{exercise}[1][]
{\begin{exmp}[#1]\begin{leftbar}}
{\end{leftbar}\end{exmp}}
\newenvironment{statement}{\begin{quote}}{\end{quote}}
\let\sumnonlimits\sum
\let\prodnonlimits\prod
\let\cupnonlimits\bigcup
\let\capnonlimits\bigcap
\renewcommand{\sum}{\sumnonlimits\limits}
\renewcommand{\prod}{\prodnonlimits\limits}
\renewcommand{\bigcup}{\cupnonlimits\limits}
\renewcommand{\bigcap}{\capnonlimits\limits}
\newenvironment{verlong}{}{}
\newenvironment{vershort}{}{}
\newenvironment{noncompile}{}{}
\excludecomment{verlong}
\includecomment{vershort}
\excludecomment{noncompile}

\ihead{Notes on the combinatorial fundamentals of algebra}
\ohead{page \thepage}
\begin{document}

\title{Notes on the combinatorial fundamentals of algebra\thanks{old title: PRIMES
2015 reading project: problems and solutions}}
\author{Darij Grinberg}
\date{September 15, 2022\\
(with minor corrections
\today
)\footnote{The numbering in this version is compatible with that in the
version of 10 January 2019 and in all intermediate versions.} }
\maketitle

\begin{abstract}
\textbf{Abstract.} This is a detailed survey -- with rigorous and
self-contained proofs -- of some of the basics of elementary combinatorics and
algebra, including the properties of finite sums, binomial coefficients,
permutations and determinants. It is entirely expository (and written to a
large extent as a repository for folklore proofs); no new results (and few, if
any, new proofs) appear.

\end{abstract}
\tableofcontents

\doclicenseThis

\section{\label{chp.intro}Introduction}

These notes are a detailed introduction to some of the basic objects of
combinatorics and algebra: finite sums, binomial coefficients, permutations
and determinants (from a combinatorial viewpoint -- no linear algebra is
presumed). To a lesser extent, modular arithmetic and recurrent integer
sequences are treated as well. The reader is assumed to be proficient in
high-school mathematics, and mature enough to understand nontrivial
mathematical proofs. Familiarity with ``contest mathematics'' is also useful.

One feature of these notes is their focus on rigorous and detailed proofs.
Indeed, so extensive are the details that a reader with experience in
mathematics will probably be able to skip whole paragraphs of proof without
losing the thread. (As a consequence of this amount of detail, the notes
contain far less material than might be expected from their length.) Rigorous
proofs mean that (with some minor exceptions) no \textquotedblleft
handwaving\textquotedblright\ is used; all relevant objects are defined in
mathematical (usually set-theoretical) language, and are manipulated in
logically well-defined ways. (In particular, some things that are commonly
taken for granted in the literature -- e.g., the fact that the sum of $n$
numbers is well-defined without specifying in what order they are being added
-- are unpacked and proven in a rigorous way.)

These notes are split into several chapters:

\begin{itemize}
\item Chapter \ref{chp.intro} collects some basic facts and notations that are
used in later chapters. This chapter is \textbf{not} meant to be read first;
it is best consulted when needed.

\item Chapter \ref{chp.ind} is an in-depth look at mathematical induction (in
various forms, including strong and two-sided induction) and several of its
applications (including basic modular arithmetic, division with remainder,
Bezout's theorem, some properties of recurrent sequences, the well-definedness
of compositions of $n$ maps and sums of $n$ numbers, and various properties thereof).

\item Chapter \ref{chp.binom} surveys binomial coefficients and their basic
properties. Unlike most texts on combinatorics, our treatment of binomial
coefficients leans to the algebraic side, relying mostly on computation and
manipulations of sums; but some basics of counting are included.

\item Chapter \ref{chp.recur} treats some more properties of Fibonacci-like
sequences, including explicit formulas (\`{a} la Binet) for two-term
recursions of the form $x_{n}=ax_{n-1}+bx_{n-2}$.

\item Chapter \ref{chp.perm} is concerned with permutations of finite sets.
The coverage is heavily influenced by the needs of the next chapter (on
determinants); thus, a great role is played by transpositions and the
inversions of a permutation.

\item Chapter \ref{chp.det} is a comprehensive introduction to determinants of
square matrices over a commutative ring\footnote{The notion of a commutative
ring is defined (and illustrated with several examples) in Section
\ref{sect.commring}, but I don't delve deeper into abstract algebra.}, from an
elementary point of view. This is probably the most unique feature of these
notes: I define determinants using Leibniz's formula (i.e., as sums over
permutations) and prove all their properties (Laplace expansion in one or
several rows; the Cauchy-Binet, Desnanot-Jacobi and Pl\"{u}cker identities;
the Vandermonde and Cauchy determinants; and several more) from this vantage
point, thus treating them as an elementary object unmoored from its
linear-algebraic origins and applications. No use is made of modules (or
vector spaces), exterior powers, eigenvalues, or of the \textquotedblleft
universal coefficients\textquotedblright\ trick\footnote{This refers to the
standard trick used for proving determinant identities (and other polynomial
identities), in which one first replaces the entries of a matrix (or, more
generally, the variables appearing in the identity) by indeterminates, then
uses the \textquotedblleft genericity\textquotedblright\ of these
indeterminates (e.g., to invert the matrix, or to divide by an expression that
could otherwise be $0$), and finally substitutes the old variables back for
the indeterminates.}. (This means that all proofs are done through
combinatorics and manipulation of sums -- a rather restrictive requirement!)
This is a conscious and (to a large extent) aesthetic choice on my part, and I
do \textbf{not} consider it the best way to learn about determinants; but I do
regard it as a road worth charting, and these notes are my attempt at doing so.
\end{itemize}

The notes include numerous exercises of varying difficulty, many of them
solved. The reader should treat exercises and theorems (and propositions,
lemmas and corollaries) as interchangeable to some extent; it is perfectly
reasonable to read the solution of an exercise, or conversely, to prove a
theorem on one's own instead of reading its proof. The reader's experience
will be the strongest determinant of their success in solving the exercises independently.

I have not meant these notes to be a textbook on any particular subject. For
one thing, their content does not map to any of the standard university
courses, but rather straddles various subjects:

\begin{itemize}
\item Much of Chapter \ref{chp.binom} (on binomial coefficients) and Chapter
\ref{chp.perm} (on permutations) is seen in a typical combinatorics class; but
my focus is more on the algebraic side and not so much on the combinatorics.

\item Chapter \ref{chp.det} studies determinants far beyond what a usual class
on linear algebra would do; but it does not include any of the other topics
that a linear algebra class usually covers (such as row reduction, vector
spaces, linear maps, eigenvectors, tensors or bilinear forms).

\item Being devoted to mathematical induction, Chapter \ref{chp.ind} appears
to cover the same ground as a typical \textquotedblleft introduction to
proofs\textquotedblright\ textbook or class (or at least one of its main
topics). In reality, however, it complements rather than competes with most
\textquotedblleft introduction to proofs\textquotedblright\ texts I have seen;
the examples I give are (with a few exceptions) nonstandard, and the focus different.

\item While the notions of rings and groups are defined in Chapter
\ref{chp.det}, I cannot claim to really be doing any abstract algebra: I am
merely working \textit{in} rings (i.e., doing computations with elements of
rings or with matrices over rings), rather than working \textit{with} rings.
Nevertheless, Chapter \ref{chp.det} might help familiarize the reader with
these concepts, facilitating proper learning of abstract algebra later on.
\end{itemize}

All in all, these notes are probably more useful as a repository of detailed
proofs than as a textbook to be read cover-to-cover. Indeed, one of my motives
in writing them was to have a reference for certain folklore results -- one in
which these results are proved elementary and without appeal to the reader's
problem-solving acumen.

These notes began as worksheets for the PRIMES reading project I have mentored
in 2015; they have since been greatly expanded with new material (some of it
originally written for my combinatorics classes, some in response to
\href{https://math.stackexchange.com/}{math.stackexchange} questions).

The notes are in flux, and probably have their share of misprints. I thank
Anya Zhang and Karthik Karnik (the two students taking part in the 2015 PRIMES
project) for finding some errors, and Christian Krattenthaler for comments.
Thanks also to the PRIMES project at MIT, which gave the impetus for the
writing of this notes; and to George Lusztig for the sponsorship of my
mentoring position in this project.

\subsection{Prerequisites}

Let me first discuss the prerequisites for a reader of these notes. At the
current moment, I assume that the reader

\begin{itemize}
\item has a good grasp on basic school-level mathematics (integers, rational
numbers, etc.);

\item has some experience with proofs (mathematical induction, proof by
contradiction, the concept of \textquotedblleft WLOG\textquotedblright, etc.)
and mathematical notation (functions, subscripts, cases, what it means for an
object to be \textquotedblleft well-defined\textquotedblright,
etc.)\footnote{A great introduction into these matters (and many others!) is
the free book \cite{LeLeMe16} by Lehman, Leighton and Meyer.
(\textbf{Practical note:} As of 2018, this book is still undergoing frequent
revisions; thus, the version I am citing below might be outdated by the time
you are reading this. I therefore suggest searching for possibly newer
versions on the internet. Unfortunately, you will also find many older
versions, often as the first google hits. Try searching for the title of the
book along with the current year to find something up-to-date.)
\par
Another introduction to proofs and mathematical workmanship is Day's
\cite{Day-proofs} (but beware that the definition of polynomials in
\cite[Chapter 5]{Day-proofs} is the wrong one for our purposes). Two others
are Hammack's \cite{Hammac15} and Doud's and Nielsen's \cite{DouNie19}. Yet
another is Newstead's \cite{Newste19} (currently a work in progress, but
promising to become one of the most interesting and sophisticated texts of
this kind). There are also several books on this subject; an especially
popular one is Velleman's \cite{Vellem06}.};

\item knows what a polynomial is (at least over $\mathbb{Z}$ and $\mathbb{Q}$)
and how polynomials differ from polynomial functions\footnote{This is used
only in a few sections and exercises, so it is not an unalienable requirement.
See Section \ref{sect.polynomials-emergency} below for a quick survey of
polynomials, and for references to sources in which precise definitions can be
found.};

\item is somewhat familiar with the summation sign ($\sum$) and the product
sign ($\prod$) and knows how to transform them (e.g., interchanging
summations, and substituting the index)\footnote{See Section
\ref{sect.sums-repetitorium} below for a quick overview of the notations that
we will need.};

\item has some familiarity with matrices (i.e., knows how to add and to
multiply them)\footnote{See, e.g., \cite[Chapter 2]{Gri-lina} or any textbook
on linear algebra for an introduction.}.
\end{itemize}

Probably a few more requirements creep in at certain points of the notes,
which I have overlooked. Some examples and remarks rely on additional
knowledge (such as analysis, graph theory, abstract algebra); however, these
can be skipped.

\subsection{Notations}

\begin{itemize}
\item In the following, we use $\mathbb{N}$ to denote the set $\left\{
0,1,2,\ldots\right\}  $. (Be warned that some other authors use the letter
$\mathbb{N}$ for $\left\{  1,2,3,\ldots\right\}  $ instead.)

\item We let $\mathbb{Q}$ denote the set of all rational numbers; we let
$\mathbb{R}$ be the set of all real numbers; we let $\mathbb{C}$ be the set of
all complex numbers\footnote{See \cite[Section 3.9]{Swanso18} or \cite[Section
I.11]{AmaEsc05} for a quick introduction to complex numbers. We will rarely
use complex numbers. Most of the time we use them, you can instead use real
numbers.}.

\item If $X$ and $Y$ are two sets, then we shall use the notation
\textquotedblleft$X\rightarrow Y,\ x\mapsto E$\textquotedblright\ (where $x$
is some symbol which has no specific meaning in the current context, and where
$E$ is some expression which usually involves $x$) for \textquotedblleft the
map from $X$ to $Y$ which sends every $x\in X$ to $E$\textquotedblright.

For example, \textquotedblleft$\mathbb{N}\rightarrow\mathbb{N},\ x\mapsto
x^{2}+x+6$\textquotedblright\ means the map from $\mathbb{N}$ to $\mathbb{N}$
which sends every $x\in\mathbb{N}$ to $x^{2}+x+6$.

For another example, \textquotedblleft$\mathbb{N}\rightarrow\mathbb{Q}%
,\ x\mapsto\dfrac{x}{1+x}$\textquotedblright\ denotes the map from
$\mathbb{N}$ to $\mathbb{Q}$ which sends every $x\in\mathbb{N}$ to $\dfrac
{x}{1+x}$.\ \ \ \ \footnote{A word of warning: Of course, the notation
\textquotedblleft$X\rightarrow Y,\ x\mapsto E$\textquotedblright\ does not
always make sense; indeed, the map that it stands for might sometimes not
exist. For instance, the notation \textquotedblleft$\mathbb{N}\rightarrow
\mathbb{Q},\ x\mapsto\dfrac{x}{1-x}$\textquotedblright\ does not actually
define a map, because the map that it is supposed to define (i.e., the map
from $\mathbb{N}$ to $\mathbb{Q}$ which sends every $x\in\mathbb{N}$ to
$\dfrac{x}{1-x}$) does not exist (since $\dfrac{x}{1-x}$ is not defined for
$x=1$). For another example, the notation \textquotedblleft$\mathbb{N}%
\rightarrow\mathbb{Z},\ x\mapsto\dfrac{x}{1+x}$\textquotedblright\ does not
define a map, because the map that it is supposed to define (i.e., the map
from $\mathbb{N}$ to $\mathbb{Z}$ which sends every $x\in\mathbb{N}$ to
$\dfrac{x}{1+x}$) does not exist (for $x=2$, we have $\dfrac{x}{1+x}=\dfrac
{2}{1+2}\notin\mathbb{Z}$, which shows that a map from $\mathbb{N}$ to
$\mathbb{Z}$ cannot send this $x$ to this $\dfrac{x}{1+x}$). Thus, when
defining a map from $X$ to $Y$ (using whatever notation), do not forget to
check that it is well-defined (i.e., that your definition specifies precisely
one image for each $x\in X$, and that these images all lie in $Y$). In many
cases, this is obvious or very easy to check (I will usually not even mention
this check), but in some cases, this is a difficult task.}

\item If $S$ is a set, then the \textit{powerset} of $S$ means the set of all
subsets of $S$. This powerset will be denoted by $\mathcal{P}\left(  S\right)
$. For example, the powerset of $\left\{  1,2\right\}  $ is $\mathcal{P}%
\left(  \left\{  1,2\right\}  \right)  =\left\{  \varnothing,\left\{
1\right\}  ,\left\{  2\right\}  ,\left\{  1,2\right\}  \right\}  $.

\item The letter $i$ will \textbf{not} denote the imaginary unit $\sqrt{-1}$
(except when we explicitly say so).
\end{itemize}

Further notations will be defined whenever they arise for the first time.

\subsection{\label{sect.jectivity}Injectivity, surjectivity, bijectivity}

In this section\footnote{a significant part of which is copied from
\cite[\S 3.21]{Gri-lina}}, we recall some basic properties of maps --
specifically, what it means for a map to be injective, surjective and
bijective. We begin by recalling basic definitions:

\begin{itemize}
\item The words \textquotedblleft map\textquotedblright, \textquotedblleft
mapping\textquotedblright, \textquotedblleft function\textquotedblright,
\textquotedblleft transformation\textquotedblright\ and \textquotedblleft
operator\textquotedblright\ are synonyms in mathematics.\footnote{That said,
mathematicians often show some nuance by using one of them and not the other.
However, we do not need to concern ourselves with this here.}

\item A map $f:X\rightarrow Y$ between two sets $X$ and $Y$ is said to be
\textit{injective} if it has the following property:

\begin{itemize}
\item If $x_{1}$ and $x_{2}$ are two elements of $X$ satisfying $f\left(
x_{1}\right)  =f\left(  x_{2}\right)  $, then $x_{1}=x_{2}$. (In words: If two
elements of $X$ are sent to one and the same element of $Y$ by $f$, then these
two elements of $X$ must have been equal in the first place. In other words:
An element of $X$ is uniquely determined by its image under $f$.)
\end{itemize}

Injective maps are often called \textquotedblleft one-to-one
maps\textquotedblright\ or \textquotedblleft injections\textquotedblright.

For example:

\begin{itemize}
\item The map $\mathbb{Z}\rightarrow\mathbb{Z},\ x\mapsto2x$ (this is the map
that sends each integer $x$ to $2x$) is injective, because if $x_{1}$ and
$x_{2}$ are two integers satisfying $2x_{1}=2x_{2}$, then $x_{1}=x_{2}$.

\item The map $\mathbb{Z}\rightarrow\mathbb{Z},\ x\mapsto x^{2}$ (this is the
map that sends each integer $x$ to $x^{2}$) is \textbf{not} injective, because
if $x_{1}$ and $x_{2}$ are two integers satisfying $x_{1}^{2}=x_{2}^{2}$, then
we do not necessarily have $x_{1}=x_{2}$. (For example, if $x_{1}=-1$ and
$x_{2}=1$, then $x_{1}^{2}=x_{2}^{2}$ but not $x_{1}=x_{2}$.)
\end{itemize}

\item A map $f:X\rightarrow Y$ between two sets $X$ and $Y$ is said to be
\textit{surjective} if it has the following property:

\begin{itemize}
\item For each $y\in Y$, there exists some $x\in X$ satisfying $f\left(
x\right)  =y$. (In words: Each element of $Y$ is an image of some element of
$X$ under $f$.)
\end{itemize}

Surjective maps are often called \textquotedblleft onto maps\textquotedblright%
\ or \textquotedblleft surjections\textquotedblright.

For example:

\begin{itemize}
\item The map $\mathbb{Z}\rightarrow\mathbb{Z},\ x\mapsto x+1$ (this is the
map that sends each integer $x$ to $x+1$) is surjective, because each integer
$y$ has some integer satisfying $x+1=y$ (namely, $x=y-1$).

\item The map $\mathbb{Z}\rightarrow\mathbb{Z},\ x\mapsto2x$ (this is the map
that sends each integer $x$ to $2x$) is \textbf{not} surjective, because not
each integer $y$ has some integer $x$ satisfying $2x=y$. (For instance, $y=1$
has no such $x$, since $y$ is odd.)

\item The map $\left\{  1,2,3,4\right\}  \rightarrow\left\{
1,2,3,4,5\right\}  ,\ x\mapsto x$ (this is the map sending each $x$ to $x$) is
\textbf{not} surjective, because not each $y\in\left\{  1,2,3,4,5\right\}  $
has some $x\in\left\{  1,2,3,4\right\}  $ satisfying $x=y$. (Namely, $y=5$ has
no such $x$.)
\end{itemize}

\item A map $f:X\rightarrow Y$ between two sets $X$ and $Y$ is said to be
\textit{bijective} if it is both injective and surjective. Bijective maps are
often called \textquotedblleft one-to-one correspondences\textquotedblright%
\ or \textquotedblleft bijections\textquotedblright.

For example:

\begin{itemize}
\item The map $\mathbb{Z}\rightarrow\mathbb{Z},\ x\mapsto x+1$ is bijective,
since it is both injective and surjective.

\item The map $\left\{  1,2,3,4\right\}  \rightarrow\left\{
1,2,3,4,5\right\}  ,\ x\mapsto x$ is \textbf{not} bijective, since it is not
surjective. (However, it is injective.)

\item The map $\mathbb{Z}\rightarrow\mathbb{N},\ x\mapsto\left\vert
x\right\vert $ is \textbf{not} bijective, since it is not injective. (However,
it is surjective.)

\item The map $\mathbb{Z}\rightarrow\mathbb{Z},\ x\mapsto x^{2}$ is
\textbf{not} bijective, since it is not injective. (It also is not surjective.)
\end{itemize}

\item If $X$ is a set, then $\operatorname*{id}\nolimits_{X}$ denotes the map
from $X$ to $X$ that sends each $x\in X$ to $x$ itself. (In words:
$\operatorname*{id}\nolimits_{X}$ denotes the map which sends each element of
$X$ to itself.) The map $\operatorname*{id}\nolimits_{X}$ is often called the
\textit{identity map on }$X$, and often denoted by $\operatorname*{id}$ (when
$X$ is clear from the context or irrelevant). The identity map
$\operatorname*{id}\nolimits_{X}$ is always bijective.

\item If $f:X\rightarrow Y$ and $g:Y\rightarrow Z$ are two maps, then the
\textit{composition} $g\circ f$ of the maps $g$ and $f$ is defined to be the
map from $X$ to $Z$ that sends each $x\in X$ to $g\left(  f\left(  x\right)
\right)  $. (In words: The composition $g\circ f$ is the map from $X$ to $Z$
that applies the map $f$ \textbf{first} and \textbf{then} applies the map
$g$.) You might find it confusing that this map is denoted by $g\circ f$
(rather than $f\circ g$), given that it proceeds by applying $f$ first and $g$
last; however, this has its reasons: It satisfies $\left(  g\circ f\right)
\left(  x\right)  =g\left(  f\left(  x\right)  \right)  $. Had we denoted it
by $f\circ g$ instead, this equality would instead become $\left(  f\circ
g\right)  \left(  x\right)  =g\left(  f\left(  x\right)  \right)  $, which
would be even more confusing.

\item If $f:X\rightarrow Y$ is a map between two sets $X$ and $Y$, then an
\textit{inverse} of $f$ means a map $g:Y\rightarrow X$ satisfying $f\circ
g=\operatorname*{id}\nolimits_{Y}$ and $g\circ f=\operatorname*{id}%
\nolimits_{X}$. (In words, the condition \textquotedblleft$f\circ
g=\operatorname*{id}\nolimits_{Y}$\textquotedblright\ means \textquotedblleft
if you start with some element $y\in Y$, then apply $g$, then apply $f$, then
you get $y$ back\textquotedblright, or equivalently \textquotedblleft the map
$f$ undoes the map $g$\textquotedblright. Similarly, the condition
\textquotedblleft$g\circ f=\operatorname*{id}\nolimits_{X}$\textquotedblright%
\ means \textquotedblleft if you start with some element $x\in X$, then apply
$f$, then apply $g$, then you get $x$ back\textquotedblright, or equivalently
\textquotedblleft the map $g$ undoes the map $f$\textquotedblright. Thus, an
inverse of $f$ means a map $g:Y\rightarrow X$ that both undoes and is undone
by $f$.)

The map $f:X\rightarrow Y$ is said to be \textit{invertible} if and only if an
inverse of $f$ exists. If an inverse of $f$ exists, then it is
unique\footnote{\textit{Proof.} Let $g_{1}$ and $g_{2}$ be two inverses of
$f$. We shall show that $g_{1}=g_{2}$.
\par
We know that $g_{1}$ is an inverse of $f$. In other words, $g_{1}$ is a map
$Y\rightarrow X$ satisfying $f\circ g_{1}=\operatorname*{id}\nolimits_{Y}$ and
$g_{1}\circ f=\operatorname*{id}\nolimits_{X}$.
\par
We know that $g_{2}$ is an inverse of $f$. In other words, $g_{2}$ is a map
$Y\rightarrow X$ satisfying $f\circ g_{2}=\operatorname*{id}\nolimits_{Y}$ and
$g_{2}\circ f=\operatorname*{id}\nolimits_{X}$.
\par
A well-known fact (known as \textit{associativity of map composition}, and
stated explicitly as Proposition \ref{prop.ind.gen-ass-maps.fgh} below) says
that if $X$, $Y$, $Z$ and $W$ are four sets, and if $c:X\rightarrow Y$,
$b:Y\rightarrow Z$ and $a:Z\rightarrow W$ are three maps, then%
\[
\left(  a\circ b\right)  \circ c=a\circ\left(  b\circ c\right)  .
\]
Applying this fact to $Y$, $X$, $Y$, $X$, $g_{1}$, $f$ and $g_{2}$ instead of
$X$, $Y$, $Z$, $W$, $c$, $b$ and $a$, we obtain $\left(  g_{2}\circ f\right)
\circ g_{1}=g_{2}\circ\left(  f\circ g_{1}\right)  $.
\par
Hence, $g_{2}\circ\left(  f\circ g_{1}\right)  =\underbrace{\left(  g_{2}\circ
f\right)  }_{=\operatorname*{id}\nolimits_{X}}\circ g_{1}=\operatorname*{id}%
\nolimits_{X}\circ g_{1}=g_{1}$. Comparing this with $g_{2}\circ
\underbrace{\left(  f\circ g_{1}\right)  }_{=\operatorname*{id}\nolimits_{Y}%
}=g_{2}\circ\operatorname*{id}\nolimits_{Y}=g_{2}$, we obtain $g_{1}=g_{2}$.
\par
Now, forget that we fixed $g_{1}$ and $g_{2}$. We thus have shown that if
$g_{1}$ and $g_{2}$ are two inverses of $f$, then $g_{1}=g_{2}$. In other
words, any two inverses of $f$ must be equal. In other words, if an inverse of
$f$ exists, then it is unique.}, and thus is called \textit{the inverse of
}$f$, and is denoted by $f^{-1}$.

For example:

\begin{itemize}
\item The map $\mathbb{Z}\rightarrow\mathbb{Z},\ x\mapsto x+1$ is invertible,
and its inverse is $\mathbb{Z}\rightarrow\mathbb{Z},\ x\mapsto x-1$.

\item The map $\mathbb{Q}\setminus\left\{  1\right\}  \rightarrow
\mathbb{Q}\setminus\left\{  0\right\}  ,\ x\mapsto\dfrac{1}{1-x}$ is
invertible, and its inverse is the map $\mathbb{Q}\setminus\left\{  0\right\}
\rightarrow\mathbb{Q}\setminus\left\{  1\right\}  ,\ x\mapsto1-\dfrac{1}{x}$.
\end{itemize}

\item If $f:X\rightarrow Y$ is a map between two sets $X$ and $Y$, then the
following notations will be used:

\begin{itemize}
\item For any subset $U$ of $X$, we let $f\left(  U\right)  $ be the subset
$\left\{  f\left(  u\right)  \ \mid\ u\in U\right\}  $ of $Y$. This set
$f\left(  U\right)  $ is called the \textit{image} of $U$ under $f$. This
should not be confused with the image $f\left(  x\right)  $ of a single
element $x\in X$ under $f$.

Note that the map $f:X\rightarrow Y$ is surjective if and only if $Y=f\left(
X\right)  $. (This is easily seen to be a restatement of the definition of
\textquotedblleft surjective\textquotedblright.)

\item For any subset $V$ of $Y$, we let $f^{-1}\left(  V\right)  $ be the
subset $\left\{  u\in X\ \mid\ f\left(  u\right)  \in V\right\}  $ of $X$.
This set $f^{-1}\left(  V\right)  $ is called the \textit{preimage} of $V$
under $f$. This should not be confused with the image $f^{-1}\left(  y\right)
$ of a single element $y\in Y$ under the inverse $f^{-1}$ of $f$ (when this
inverse exists).

(Note that in general, $f\left(  f^{-1}\left(  V\right)  \right)  \neq V$ and
$f^{-1}\left(  f\left(  U\right)  \right)  \neq U$. However, $f\left(
f^{-1}\left(  V\right)  \right)  \subseteq V$ and $U\subseteq f^{-1}\left(
f\left(  U\right)  \right)  $.)

\item For any subset $U$ of $X$, we let $f\mid_{U}$ be the map from $U$ to $Y$
which sends each $u\in U$ to $f\left(  u\right)  \in Y$. This map $f\mid_{U}$
is called the \textit{restriction} of $f$ to the subset $U$.
\end{itemize}
\end{itemize}

The following facts are fundamental:

\begin{theorem}
\label{thm.bijective=invertible}A map $f:X\rightarrow Y$ is invertible if and
only if it is bijective.
\end{theorem}

\begin{theorem}
\label{thm.bijection=eqsize}Let $U$ and $V$ be two finite sets. Then,
$\left\vert U\right\vert =\left\vert V\right\vert $ if and only if there
exists a bijective map $f:U\rightarrow V$.
\end{theorem}

Theorem \ref{thm.bijection=eqsize} holds even if the sets $U$ and $V$ are
infinite, but to make sense of this we would need to define the size of an
infinite set, which is a much subtler issue than the size of a finite set. We
will only need Theorem \ref{thm.bijection=eqsize} for finite sets.

Let us state some more well-known and basic properties of maps between finite sets:

\begin{lemma}
\label{lem.jectivity.pigeon0}Let $U$ and $V$ be two finite sets. Let
$f:U\rightarrow V$ be a map.

\textbf{(a)} We have $\left\vert f\left(  S\right)  \right\vert \leq\left\vert
S\right\vert $ for each subset $S$ of $U$.

\textbf{(b)} If $\left\vert f\left(  U\right)  \right\vert \geq\left\vert
U\right\vert $, then the map $f$ is injective.

\textbf{(c)} If $f$ is injective, then $\left\vert f\left(  S\right)
\right\vert =\left\vert S\right\vert $ for each subset $S$ of $U$.
\end{lemma}

\begin{lemma}
\label{lem.jectivity.pigeon-surj}Let $U$ and $V$ be two finite sets such that
$\left\vert U\right\vert \leq\left\vert V\right\vert $. Let $f:U\rightarrow V$
be a map. Then, we have the following logical equivalence:%
\[
\left(  f\text{ is surjective}\right)  \ \Longleftrightarrow\ \left(  f\text{
is bijective}\right)  .
\]

\end{lemma}

\begin{lemma}
\label{lem.jectivity.pigeon-inj}Let $U$ and $V$ be two finite sets such that
$\left\vert U\right\vert \geq\left\vert V\right\vert $. Let $f:U\rightarrow V$
be a map. Then, we have the following logical equivalence:%
\[
\left(  f\text{ is injective}\right)  \ \Longleftrightarrow\ \left(  f\text{
is bijective}\right)  .
\]

\end{lemma}

\begin{exercise}
\label{exe.jectivity.pigeons}Prove Lemma \ref{lem.jectivity.pigeon0}, Lemma
\ref{lem.jectivity.pigeon-surj} and Lemma \ref{lem.jectivity.pigeon-inj}.
\end{exercise}

Let us make one additional observation about maps:

\begin{remark}
Composition of maps is associative: If $X$, $Y$, $Z$ and $W$ are three sets,
and if $c:X\rightarrow Y$, $b:Y\rightarrow Z$ and $a:Z\rightarrow W$ are three
maps, then $\left(  a\circ b\right)  \circ c=a\circ\left(  b\circ c\right)  $.
(This shall be proven in Proposition \ref{prop.ind.gen-ass-maps.fgh} below.)

In Section \ref{sect.ind.gen-ass}, we shall prove a more general fact: If
$X_{1},X_{2},\ldots,X_{k+1}$ are $k+1$ sets for some $k\in\mathbb{N}$, and if
$f_{i}:X_{i}\rightarrow X_{i+1}$ is a map for each $i\in\left\{
1,2,\ldots,k\right\}  $, then the composition $f_{k}\circ f_{k-1}\circ
\cdots\circ f_{1}$ of all $k$ maps $f_{1},f_{2},\ldots,f_{k}$ is a
well-defined map from $X_{1}$ to $X_{k+1}$, which sends each element $x\in
X_{1}$ to $f_{k}\left(  f_{k-1}\left(  f_{k-2}\left(  \cdots\left(
f_{2}\left(  f_{1}\left(  x\right)  \right)  \right)  \cdots\right)  \right)
\right)  $ (in other words, which transforms each element $x\in X_{1}$ by
first applying $f_{1}$, then applying $f_{2}$, then applying $f_{3}$, and so
on); this composition $f_{k}\circ f_{k-1}\circ\cdots\circ f_{1}$ can also be
written as $f_{k}\circ\left(  f_{k-1}\circ\left(  f_{k-2}\circ\left(
\cdots\circ\left(  f_{2}\circ f_{1}\right)  \cdots\right)  \right)  \right)  $
or as $\left(  \left(  \left(  \cdots\left(  f_{k}\circ f_{k-1}\right)
\circ\cdots\right)  \circ f_{3}\right)  \circ f_{2}\right)  \circ f_{1}$. An
important particular case is when $k=0$; in this case, $f_{k}\circ
f_{k-1}\circ\cdots\circ f_{1}$ is a composition of $0$ maps. It is defined to
be $\operatorname*{id}\nolimits_{X_{1}}$ (the identity map of the set $X_{1}%
$), and it is called the \textquotedblleft empty composition of maps
$X_{1}\rightarrow X_{1}$\textquotedblright. (The logic behind this definition
is that the composition $f_{k}\circ f_{k-1}\circ\cdots\circ f_{1}$ should
transform each element $x\in X_{1}$ by first applying $f_{1}$, then applying
$f_{2}$, then applying $f_{3}$, and so on; however, for $k=0$, there are no
maps to apply, and so $x$ just remains unchanged.)
\end{remark}

\subsection{\label{sect.sums-repetitorium}Sums and products: a synopsis}

In this section, I will recall the definitions of the $\sum$ and $\prod$ signs
and collect some of their basic properties (without proofs). When I say
\textquotedblleft recall\textquotedblright, I am implying that the reader has
at least some prior acquaintance (and, ideally, experience) with these signs;
for a first introduction, this section is probably too brief and too abstract.
Ideally, you should use this section to familiarize yourself with my
(sometimes idiosyncratic) notations.

Throughout Section \ref{sect.sums-repetitorium}, we let $\mathbb{A}$ be one of
the sets $\mathbb{N}$, $\mathbb{Z}$, $\mathbb{Q}$, $\mathbb{R}$ and
$\mathbb{C}$.

\subsubsection{Definition of $\sum$}

Let us first define the $\sum$ sign. There are actually several (slightly
different, but still closely related) notations involving the $\sum$ sign; let
us define the most important of them:

\begin{itemize}
\item If $S$ is a finite set, and if $a_{s}$ is an element of $\mathbb{A}$ for
each $s\in S$, then $\sum_{s\in S}a_{s}$ denotes the sum of all of these
elements $a_{s}$. Formally, this sum is defined by recursion on $\left\vert
S\right\vert $, as follows:

\begin{itemize}
\item If $\left\vert S\right\vert =0$, then $\sum_{s\in S}a_{s}$ is defined to
be $0$.

\item Let $n\in\mathbb{N}$. Assume that we have defined $\sum_{s\in S}a_{s}$
for every finite set $S$ with $\left\vert S\right\vert =n$ (and every choice
of elements $a_{s}$ of $\mathbb{A}$). Now, if $S$ is a finite set with
$\left\vert S\right\vert =n+1$ (and if $a_{s}\in\mathbb{A}$ are chosen for all
$s\in S$), then $\sum_{s\in S}a_{s}$ is defined by picking any $t\in
S$\ \ \ \ \footnote{This is possible, because $S$ is nonempty (in fact,
$\left\vert S\right\vert =n+1>n\geq0$).} and setting%
\begin{equation}
\sum_{s\in S}a_{s}=a_{t}+\sum_{s\in S\setminus\left\{  t\right\}  }a_{s}.
\label{eq.sum.def.1}%
\end{equation}
It is not immediately clear why this definition is legitimate: The right hand
side of (\ref{eq.sum.def.1}) is defined using a choice of $t$, but we want our
value of $\sum_{s\in S}a_{s}$ to depend only on $S$ and on the $a_{s}$ (not on
some arbitrarily chosen $t\in S$). However, it is possible to prove that the
right hand side of (\ref{eq.sum.def.1}) is actually independent of $t$ (that
is, any two choices of $t$ will lead to the same result). See Section
\ref{sect.ind.gen-com} below (and Theorem \ref{thm.ind.gen-com.wd}
\textbf{(a)} in particular) for the proof of this fact.
\end{itemize}

\textbf{Examples:}

\begin{itemize}
\item If $S=\left\{  4,7,9\right\}  $ and $a_{s}=\dfrac{1}{s^{2}}$ for every
$s\in S$, then $\sum_{s\in S}a_{s}=a_{4}+a_{7}+a_{9}=\dfrac{1}{4^{2}}%
+\dfrac{1}{7^{2}}+\dfrac{1}{9^{2}}=\dfrac{6049}{63504}$.

\item If $S=\left\{  1,2,\ldots,n\right\}  $ (for some $n\in\mathbb{N}$) and
$a_{s}=s^{2}$ for every $s\in S$, then $\sum_{s\in S}a_{s}=\sum_{s\in S}%
s^{2}=1^{2}+2^{2}+\cdots+n^{2}$. (There is a formula saying that the right
hand side of this equality is $\dfrac{1}{6}n\left(  2n+1\right)  \left(
n+1\right)  $.)

\item If $S=\varnothing$, then $\sum_{s\in S}a_{s}=0$ (since $\left\vert
S\right\vert =0$).
\end{itemize}

\textbf{Remarks:}

\begin{itemize}
\item The sum $\sum_{s\in S}a_{s}$ is usually pronounced \textquotedblleft sum
of the $a_{s}$ over all $s\in S$\textquotedblright\ or \textquotedblleft sum
of the $a_{s}$ with $s$ ranging over $S$\textquotedblright\ or
\textquotedblleft sum of the $a_{s}$ with $s$ running through all elements of
$S$\textquotedblright. The letter \textquotedblleft$s$\textquotedblright\ in
the sum is called the \textquotedblleft summation index\textquotedblright%
\footnote{The plural of the word \textquotedblleft index\textquotedblright%
\ here is \textquotedblleft indices\textquotedblright, not \textquotedblleft
indexes\textquotedblright.}, and its exact choice is immaterial (for example,
you can rewrite $\sum_{s\in S}a_{s}$ as $\sum_{t\in S}a_{t}$ or as $\sum
_{\Phi\in S}a_{\Phi}$ or as $\sum_{\spadesuit\in S}a_{\spadesuit}$), as long
as it does not already have a different meaning outside of the sum\footnote{If
it already has a different meaning, then it must not be used as a summation
index! For example, you must not write \textquotedblleft every $n\in
\mathbb{N}$ satisfies $\sum_{n\in\left\{  0,1,\ldots,n\right\}  }%
n=\dfrac{n\left(  n+1\right)  }{2}$\textquotedblright, because here the
summation index $n$ clashes with a different meaning of the letter $n$.}.
(Ultimately, a summation index is the same kind of placeholder variable as the
\textquotedblleft$s$\textquotedblright\ in the statement \textquotedblleft for
all $s\in S$, we have $a_{s}+2a_{s}=3a_{s}$\textquotedblright, or as a loop
variable in a for-loop in programming.) The sign $\sum$ itself is called
\textquotedblleft the summation sign\textquotedblright\ or \textquotedblleft
the $\sum$ sign\textquotedblright. The numbers $a_{s}$ are called the
\textit{addends} (or \textit{summands}) of the sum $\sum_{s\in S}a_{s}$. More
precisely, for any given $t\in S$, we can refer to the number $a_{t}$ as the
\textquotedblleft addend corresponding to the index $t$\textquotedblright\ (or
as the \textquotedblleft addend for $s=t$\textquotedblright, or as the
\textquotedblleft addend for $t$\textquotedblright) of the sum $\sum_{s\in
S}a_{s}$.

\item When the set $S$ is empty, the sum $\sum_{s\in S}a_{s}$ is called an
\textit{empty sum}. Our definition implies that any empty sum is $0$. This
convention is used throughout mathematics, except in rare occasions where a
slightly subtler version of it is used\footnote{Do not worry about this
subtler version for the time being. If you really want to know what it is: Our
above definition is tailored to the cases when the $a_{s}$ are numbers (i.e.,
elements of one of the sets $\mathbb{N}$, $\mathbb{Z}$, $\mathbb{Q}$,
$\mathbb{R}$ and $\mathbb{C}$). In more advanced settings, one tends to take
sums of the form $\sum_{s\in S}a_{s}$ where the $a_{s}$ are not numbers but
(for example) elements of a commutative ring $\mathbb{K}$. (See Definition
\ref{def.commring} for the definition of a commutative ring.) In such cases,
one wants the sum $\sum_{s\in S}a_{s}$ for an empty set $S$ to be not the
integer $0$, but the zero of the commutative ring $\mathbb{K}$ (which is
sometimes distinct from the integer $0$). This has the slightly confusing
consequence that the meaning of the sum $\sum_{s\in S}a_{s}$ for an empty set
$S$ depends on what ring $\mathbb{K}$ the $a_{s}$ belong to, even if (for an
empty set $S$) there are no $a_{s}$ to begin with! But in practice, the choice
of $\mathbb{K}$ is always clear from context, so this is not ambiguous.
\par
A similar caveat applies to the other versions of the $\sum$ sign, as well as
to the $\prod$ sign defined further below; I shall not elaborate on it
further.}. Ignore anyone who tells you that empty sums are undefined!

\item The summation index does not always have to be a single letter. For
instance, if $S$ is a set of pairs, then we can write $\sum_{\left(
x,y\right)  \in S}a_{\left(  x,y\right)  }$ (meaning the same as $\sum_{s\in
S}a_{s}$). Here is an example of this notation:%
\[
\sum_{\left(  x,y\right)  \in\left\{  1,2,3\right\}  ^{2}}\ \ \dfrac{x}%
{y}=\dfrac{1}{1}+\dfrac{1}{2}+\dfrac{1}{3}+\dfrac{2}{1}+\dfrac{2}{2}+\dfrac
{2}{3}+\dfrac{3}{1}+\dfrac{3}{2}+\dfrac{3}{3}%
\]
(here, we are using the notation $\sum_{\left(  x,y\right)  \in S}a_{\left(
x,y\right)  }$ with $S=\left\{  1,2,3\right\}  ^{2}$ and $a_{\left(
x,y\right)  }=\dfrac{x}{y}$). Note that we could not have rewritten this sum
in the form $\sum_{s\in S}a_{s}$ with a single-letter variable $s$ without
introducing an extra notation such as $a_{\left(  x,y\right)  }$ for the
quotients $\dfrac{x}{y}$.

\item Mathematicians don't seem to have reached an agreement on the operator
precedence of the $\sum$ sign. By this I mean the following question: Does
$\sum_{s\in S}a_{s}+b$ (where $b$ is some other element of $\mathbb{A}$) mean
$\sum_{s\in S}\left(  a_{s}+b\right)  $ or $\left(  \sum_{s\in S}a_{s}\right)
+b$ ? In my experience, the second interpretation (i.e., reading it as
$\left(  \sum_{s\in S}a_{s}\right)  +b$) is more widespread, and this is the
interpretation that I will follow. Nevertheless, be on the watch for possible
misunderstandings, as someone might be using the first interpretation when you
expect it the least!\footnote{This is similar to the notorious disagreement
about whether $a/bc$ means $\left(  a/b\right)  \cdot c$ or $a/\left(
bc\right)  $.}

However, the situation is different for products and nested sums. For
instance, the expression $\sum_{s\in S}ba_{s}c$ is understood to mean
$\sum_{s\in S}\left(  ba_{s}c\right)  $, and a nested sum like $\sum_{s\in
S}\sum_{t\in T}a_{s,t}$ (where $S$ and $T$ are two sets, and where $a_{s,t}$
is an element of $\mathbb{A}$ for each pair $\left(  s,t\right)  \in S\times
T$) is to be read as $\sum_{s\in S}\left(  \sum_{t\in T}a_{s,t}\right)  $.

\item Speaking of nested sums: they mean exactly what they seem to mean. For
instance, $\sum_{s\in S}\sum_{t\in T}a_{s,t}$ is what you get if you compute
the sum $\sum_{t\in T}a_{s,t}$ for each $s\in S$, and then sum up all of these
sums together. In a nested sum $\sum_{s\in S}\sum_{t\in T}a_{s,t}$, the first
summation sign ($\sum_{s\in S}$) is called the \textquotedblleft outer
summation\textquotedblright, and the second summation sign ($\sum_{t\in T}$)
is called the \textquotedblleft inner summation\textquotedblright.

\item An expression of the form \textquotedblleft$\sum_{s\in S}a_{s}%
$\textquotedblright\ (where $S$ is a finite set) is called a \textit{finite
sum}.

\item We have required the set $S$ to be finite when defining $\sum_{s\in
S}a_{s}$. Of course, this requirement was necessary for our definition, and
there is no way to make sense of infinite sums such as $\sum_{s\in\mathbb{Z}%
}s^{2}$. However, \textbf{some} infinite sums can be made sense of. The
simplest case is when the set $S$ might be infinite, but only finitely many
among the $a_{s}$ are nonzero. In this case, we can define $\sum_{s\in S}%
a_{s}$ simply by discarding the zero addends and summing the finitely many
remaining addends. Other situations in which infinite sums make sense appear
in analysis and in topological algebra (e.g., power series).

\item The sum $\sum_{s\in S}a_{s}$ always belongs to $\mathbb{A}%
$.\ \ \ \ \footnote{Recall that we have assumed $\mathbb{A}$ to be one of the
sets $\mathbb{N}$, $\mathbb{Z}$, $\mathbb{Q}$, $\mathbb{R}$ and $\mathbb{C}$,
and that we have assumed the $a_{s}$ to belong to $\mathbb{A}$.} For instance,
a sum of elements of $\mathbb{N}$ belongs to $\mathbb{N}$; a sum of elements
of $\mathbb{R}$ belongs to $\mathbb{R}$, and so on.
\end{itemize}

\item A slightly more complicated version of the summation sign is the
following: Let $S$ be a finite set, and let $\mathcal{A}\left(  s\right)  $ be
a logical statement defined for every $s\in S$\ \ \ \ \footnote{Formally
speaking, this means that $\mathcal{A}$ is a map from $S$ to the set of all
logical statements. Such a map is called a \textit{predicate}.}. For example,
$S$ can be $\left\{  1,2,3,4\right\}  $, and $\mathcal{A}\left(  s\right)  $
can be the statement \textquotedblleft$s$ is even\textquotedblright. For each
$s\in S$ satisfying $\mathcal{A}\left(  s\right)  $, let $a_{s}$ be an element
of $\mathbb{A}$. Then, the sum $\sum_{\substack{s\in S;\\\mathcal{A}\left(
s\right)  }}a_{s}$ is defined by%
\[
\sum_{\substack{s\in S;\\\mathcal{A}\left(  s\right)  }}a_{s}=\sum
_{s\in\left\{  t\in S\ \mid\ \mathcal{A}\left(  t\right)  \right\}  }a_{s}.
\]
In other words, $\sum_{\substack{s\in S;\\\mathcal{A}\left(  s\right)  }%
}a_{s}$ is the sum of the $a_{s}$ for all $s\in S$ which satisfy
$\mathcal{A}\left(  s\right)  $.

\textbf{Examples:}

\begin{itemize}
\item If $S=\left\{  1,2,3,4,5\right\}  $, then $\sum_{\substack{s\in
S;\\s\text{ is even}}}a_{s}=a_{2}+a_{4}$. (Of course, $\sum_{\substack{s\in
S;\\s\text{ is even}}}a_{s}$ is $\sum_{\substack{s\in S;\\\mathcal{A}\left(
s\right)  }}a_{s}$ when $\mathcal{A}\left(  s\right)  $ is defined to be the
statement \textquotedblleft$s$ is even\textquotedblright.)

\item If $S=\left\{  1,2,\ldots,n\right\}  $ (for some $n\in\mathbb{N}$) and
$a_{s}=s^{2}$ for every $s\in S$, then $\sum_{\substack{s\in S;\\s\text{ is
even}}}a_{s}=a_{2}+a_{4}+\cdots+a_{k}$, where $k$ is the largest even number
among $1,2,\ldots,n$ (that is, $k=n$ if $n$ is even, and $k=n-1$ otherwise).
\end{itemize}

\textbf{Remarks:}

\begin{itemize}
\item The sum $\sum_{\substack{s\in S;\\\mathcal{A}\left(  s\right)  }}a_{s}$
is usually pronounced \textquotedblleft sum of the $a_{s}$ over all $s\in S$
satisfying $\mathcal{A}\left(  s\right)  $\textquotedblright. The semicolon
after \textquotedblleft$s\in S$\textquotedblright\ is often omitted or
replaced by a colon or a comma. Many authors often omit the \textquotedblleft%
$s\in S$\textquotedblright\ part (so they simply write $\sum_{\mathcal{A}%
\left(  s\right)  }a_{s}$) when it is clear enough what the $S$ is. (For
instance, they would write $\sum_{1\leq s\leq5}s^{2}$ instead of
$\sum_{\substack{s\in\mathbb{N};\\1\leq s\leq5}}s^{2}$.)

\item The set $S$ needs not be finite in order for $\sum_{\substack{s\in
S;\\\mathcal{A}\left(  s\right)  }}a_{s}$ to be defined; it suffices that the
set $\left\{  t\in S\ \mid\ \mathcal{A}\left(  t\right)  \right\}  $ be finite
(i.e., that only finitely many $s\in S$ satisfy $\mathcal{A}\left(  s\right)
$).

\item The sum $\sum_{\substack{s\in S;\\\mathcal{A}\left(  s\right)  }}a_{s}$
is said to be \textit{empty} whenever the set $\left\{  t\in S\ \mid
\ \mathcal{A}\left(  t\right)  \right\}  $ is empty (i.e., whenever no $s\in
S$ satisfies $\mathcal{A}\left(  s\right)  $).
\end{itemize}

\item Finally, here is the simplest version of the summation sign: Let $u$ and
$v$ be two integers. We agree to understand the set $\left\{  u,u+1,\ldots
,v\right\}  $ to be empty when $u>v$. Let $a_{s}$ be an element of
$\mathbb{A}$ for each $s\in\left\{  u,u+1,\ldots,v\right\}  $. Then,
$\sum_{s=u}^{v}a_{s}$ is defined by%
\[
\sum_{s=u}^{v}a_{s}=\sum_{s\in\left\{  u,u+1,\ldots,v\right\}  }a_{s}.
\]

\textbf{Examples:}

\begin{itemize}
\item We have $\sum_{s=3}^{8}\dfrac{1}{s}=\sum_{s\in\left\{  3,4,\ldots
,8\right\}  }\dfrac{1}{s}=\dfrac{1}{3}+\dfrac{1}{4}+\dfrac{1}{5}+\dfrac{1}%
{6}+\dfrac{1}{7}+\dfrac{1}{8}=\dfrac{341}{280}$.

\item We have $\sum_{s=3}^{3}\dfrac{1}{s}=\sum_{s\in\left\{  3\right\}
}\dfrac{1}{s}=\dfrac{1}{3}$.

\item We have $\sum_{s=3}^{2}\dfrac{1}{s}=\sum_{s\in\varnothing}\dfrac{1}%
{s}=0$.
\end{itemize}

\textbf{Remarks:}

\begin{itemize}
\item The sum $\sum_{s=u}^{v}a_{s}$ is usually pronounced \textquotedblleft
sum of the $a_{s}$ for all $s$ from $u$ to $v$ (inclusive)\textquotedblright.
It is often written $a_{u}+a_{u+1}+\cdots+a_{v}$, but this latter notation has
its drawbacks: In order to understand an expression like $a_{u}+a_{u+1}%
+\cdots+a_{v}$, one needs to correctly guess the pattern (which can be
unintuitive when the $a_{s}$ themselves are complicated: for example, it takes
a while to find the \textquotedblleft moving parts\textquotedblright\ in the
expression $\dfrac{2\cdot7}{3+2}+\dfrac{3\cdot7}{3+3}+\cdots+\dfrac{7\cdot
7}{3+7}$, whereas the notation $\sum_{s=2}^{7}\dfrac{s\cdot7}{3+s}$ for the
same sum is perfectly clear).

\item In the sum $\sum_{s=u}^{v}a_{s}$, the integer $u$ is called the
\textit{lower limit} (of the sum), whereas the integer $v$ is called the
\textit{upper limit} (of the sum). The sum is said to \textit{start} (or
\textit{begin}) at $u$ and \textit{end} at $v$.

\item The sum $\sum_{s=u}^{v}a_{s}$ is said to be \textit{empty} whenever
$u>v$. In other words, a sum of the form $\sum_{s=u}^{v}a_{s}$ is empty
whenever it \textquotedblleft ends before it has begun\textquotedblright.
However, a sum which \textquotedblleft ends right after it
begins\textquotedblright\ (i.e., a sum $\sum_{s=u}^{v}a_{s}$ with $u=v$) is
not empty; it just has one addend only. (This is unlike integrals, which are
$0$ whenever their lower and upper limit are equal.)

\item Let me stress once again that a sum $\sum_{s=u}^{v}a_{s}$ with $u>v$ is
empty and equals $0$. It does not matter how much greater $u$ is than $v$. So,
for example, $\sum_{s=1}^{-5}s=0$. The fact that the upper bound ($-5$) is
much smaller than the lower bound ($1$) does not mean that you have to
subtract rather than add.
\end{itemize}
\end{itemize}

Thus we have introduced the main three forms of the summation sign. Some mild
variations on them appear in the literature (e.g., there is a slightly awkward
notation $\sum_{\substack{s=u;\\\mathcal{A}\left(  s\right)  }}^{v}a_{s}$ for
$\sum_{\substack{s\in\left\{  u,u+1,\ldots,v\right\}  ;\\\mathcal{A}\left(
s\right)  }}a_{s}$).

\subsubsection{Properties of $\sum$}

Let me now show some basic properties of summation signs that are important in
making them useful:

\begin{itemize}
\item \underline{\textbf{Splitting-off:}} Let $S$ be a finite set. Let $t\in
S$. Let $a_{s}$ be an element of $\mathbb{A}$ for each $s\in S$. Then,%
\begin{equation}
\sum_{s\in S}a_{s}=a_{t}+\sum_{s\in S\setminus\left\{  t\right\}  }a_{s}.
\label{eq.sum.split-off}%
\end{equation}
(This is precisely the equality (\ref{eq.sum.def.1}) (applied to $n=\left\vert
S\setminus\left\{  t\right\}  \right\vert $), because $\left\vert S\right\vert
=\left\vert S\setminus\left\{  t\right\}  \right\vert +1$.) This formula
(\ref{eq.sum.split-off}) allows us to \textquotedblleft split
off\textquotedblright\ an addend from a sum.

\textbf{Example:} If $n\in\mathbb{N}$, then
\[
\sum_{s\in\left\{  1,2,\ldots,n+1\right\}  }a_{s}=a_{n+1}+\sum_{s\in\left\{
1,2,\ldots,n\right\}  }a_{s}%
\]
(by (\ref{eq.sum.split-off}), applied to $S=\left\{  1,2,\ldots,n+1\right\}  $
and $t=n+1$), but also%
\[
\sum_{s\in\left\{  1,2,\ldots,n+1\right\}  }a_{s}=a_{1}+\sum_{s\in\left\{
2,3,\ldots,n+1\right\}  }a_{s}%
\]
(by (\ref{eq.sum.split-off}), applied to $S=\left\{  1,2,\ldots,n+1\right\}  $
and $t=1$).

\item \underline{\textbf{Splitting:}} Let $S$ be a finite set. Let $X$ and $Y$
be two subsets of $S$ such that $X\cap Y=\varnothing$ and $X\cup Y=S$.
(Equivalently, $X$ and $Y$ are two subsets of $S$ such that each element of
$S$ lies in \textbf{exactly} one of $X$ and $Y$.) Let $a_{s}$ be an element of
$\mathbb{A}$ for each $s\in S$. Then,%
\begin{equation}
\sum_{s\in S}a_{s}=\sum_{s\in X}a_{s}+\sum_{s\in Y}a_{s}. \label{eq.sum.split}%
\end{equation}
(Here, as we explained, $\sum_{s\in X}a_{s}+\sum_{s\in Y}a_{s}$ stands for
$\left(  \sum_{s\in X}a_{s}\right)  +\left(  \sum_{s\in Y}a_{s}\right)  $.)
The idea behind (\ref{eq.sum.split}) is that if we want to add a bunch of
numbers (the $a_{s}$ for $s\in S$), we can proceed by splitting it into two
\textquotedblleft sub-bunches\textquotedblright\ (one \textquotedblleft
sub-bunch\textquotedblright\ consisting of the $a_{s}$ for $s\in X$, and the
other consisting of the $a_{s}$ for $s\in Y$), then take the sum of each of
these two sub-bunches, and finally add together the two sums. For a rigorous
proof of (\ref{eq.sum.split}), see Theorem \ref{thm.ind.gen-com.split2} below.

\textbf{Examples:}

\begin{itemize}
\item If $n\in\mathbb{N}$, then%
\[
\sum_{s\in\left\{  1,2,\ldots,2n\right\}  }a_{s}=\sum_{s\in\left\{
1,3,\ldots,2n-1\right\}  }a_{s}+\sum_{s\in\left\{  2,4,\ldots,2n\right\}
}a_{s}%
\]
(by (\ref{eq.sum.split}), applied to $S=\left\{  1,2,\ldots,2n\right\}  $,
$X=\left\{  1,3,\ldots,2n-1\right\}  $ and $Y=\left\{  2,4,\ldots,2n\right\}
$).

\item If $n\in\mathbb{N}$ and $m\in\mathbb{N}$, then%
\[
\sum_{s\in\left\{  -m,-m+1,\ldots,n\right\}  }a_{s}=\sum_{s\in\left\{
-m,-m+1,\ldots,0\right\}  }a_{s}+\sum_{s\in\left\{  1,2,\ldots,n\right\}
}a_{s}%
\]
(by (\ref{eq.sum.split}), applied to $S=\left\{  -m,-m+1,\ldots,n\right\}  $,
$X=\left\{  -m,-m+1,\ldots,0\right\}  $ and $Y=\left\{  1,2,\ldots,n\right\}
$).

\item If $u$, $v$ and $w$ are three integers such that $u-1\leq v\leq w$, and
if $a_{s}$ is an element of $\mathbb{A}$ for each $s\in\left\{  u,u+1,\ldots
,w\right\}  $, then%
\begin{equation}
\sum_{s=u}^{w}a_{s}=\sum_{s=u}^{v}a_{s}+\sum_{s=v+1}^{w}a_{s}.
\label{eq.sum.split.uvw}%
\end{equation}
This follows from (\ref{eq.sum.split}), applied to $S=\left\{  u,u+1,\ldots
,w\right\}  $, $X=\left\{  u,u+1,\ldots,v\right\}  $ and $Y=\left\{
v+1,v+2,\ldots,w\right\}  $. Notice that the requirement $u-1\leq v\leq w$ is
important; otherwise, the $X\cap Y=\varnothing$ and $X\cup Y=S$ condition
would not hold!
\end{itemize}

\item \underline{\textbf{Splitting using a predicate:}} Let $S$ be a finite
set. Let $\mathcal{A}\left(  s\right)  $ be a logical statement for each $s\in
S$. Let $a_{s}$ be an element of $\mathbb{A}$ for each $s\in S$. Then,%
\begin{equation}
\sum_{s\in S}a_{s}=\sum_{\substack{s\in S;\\\mathcal{A}\left(  s\right)
}}a_{s}+\sum_{\substack{s\in S;\\\text{not }\mathcal{A}\left(  s\right)
}}a_{s} \label{eq.sum.split.pred}%
\end{equation}
(where \textquotedblleft not $\mathcal{A}\left(  s\right)  $\textquotedblright%
\ means the negation of $\mathcal{A}\left(  s\right)  $). This simply follows
from (\ref{eq.sum.split}), applied to $X=\left\{  s\in S\ \mid\ \mathcal{A}%
\left(  s\right)  \right\}  $ and $Y=\left\{  s\in S\ \mid\ \text{not
}\mathcal{A}\left(  s\right)  \right\}  $.

\textbf{Example:} If $S\subseteq\mathbb{Z}$, then%
\[
\sum_{s\in S}a_{s}=\sum_{\substack{s\in S;\\s\text{ is even}}}a_{s}%
+\sum_{\substack{s\in S;\\s\text{ is odd}}}a_{s}%
\]
(because \textquotedblleft$s$ is odd\textquotedblright\ is the negation of
\textquotedblleft$s$ is even\textquotedblright).

\item \underline{\textbf{Summing equal values:}} Let $S$ be a finite set. Let
$a$ be an element of $\mathbb{A}$. Then,%
\begin{equation}
\sum_{s\in S}a=\left\vert S\right\vert \cdot a. \label{eq.sum.equal}%
\end{equation}
\footnote{This is easy to prove by induction on $\left\vert S\right\vert $.}
In other words, if all addends of a sum are equal to one and the same element
$a$, then the sum is just the number of its addends times $a$. In particular,%
\[
\sum_{s\in S}1=\left\vert S\right\vert \cdot1=\left\vert S\right\vert .
\]

\item \underline{\textbf{Splitting an addend:}} Let $S$ be a finite set. For
every $s\in S$, let $a_{s}$ and $b_{s}$ be elements of $\mathbb{A}$. Then,%
\begin{equation}
\sum_{s\in S}\left(  a_{s}+b_{s}\right)  =\sum_{s\in S}a_{s}+\sum_{s\in
S}b_{s}. \label{eq.sum.linear1}%
\end{equation}
For a rigorous proof of this equality, see Theorem
\ref{thm.ind.gen-com.sum(a+b)} below.

\textbf{Remark:} Of course, similar rules hold for other forms of summations:
If $\mathcal{A}\left(  s\right)  $ is a logical statement for each $s\in S$,
then%
\[
\sum_{\substack{s\in S;\\\mathcal{A}\left(  s\right)  }}\left(  a_{s}%
+b_{s}\right)  =\sum_{\substack{s\in S;\\\mathcal{A}\left(  s\right)  }%
}a_{s}+\sum_{\substack{s\in S;\\\mathcal{A}\left(  s\right)  }}b_{s}.
\]
If $u$ and $v$ are two integers, then%
\begin{equation}
\sum_{s=u}^{v}\left(  a_{s}+b_{s}\right)  =\sum_{s=u}^{v}a_{s}+\sum_{s=u}%
^{v}b_{s}. \label{eq.sum.linear1.c}%
\end{equation}

\item \underline{\textbf{Factoring out:}} Let $S$ be a finite set. For every
$s\in S$, let $a_{s}$ be an element of $\mathbb{A}$. Also, let $\lambda$ be an
element of $\mathbb{A}$. Then,%
\begin{equation}
\sum_{s\in S}\lambda a_{s}=\lambda\sum_{s\in S}a_{s}. \label{eq.sum.linear2}%
\end{equation}
For a rigorous proof of this equality, see Theorem
\ref{thm.ind.gen-com.sum(la)} below.

Again, similar rules hold for the other types of summation sign.

\textbf{Remark:} Applying (\ref{eq.sum.linear2}) to $\lambda=-1$, we obtain%
\[
\sum_{s\in S}\left(  -a_{s}\right)  =-\sum_{s\in S}a_{s}.
\]

\item \underline{\textbf{Zeroes sum to zero:}} Let $S$ be a finite set. Then,%
\begin{equation}
\sum_{s\in S}0=0. \label{eq.sum.sum0}%
\end{equation}
That is, any sum of zeroes is zero.

For a rigorous proof of this equality, see Theorem
\ref{thm.ind.gen-com.sum(0)} below.

\textbf{Remark:} This applies even to infinite sums! Do not be fooled by the
infiniteness of a sum: There are no reasonable situations where an infinite
sum of zeroes is defined to be anything other than zero. The infinity does not
\textquotedblleft compensate\textquotedblright\ for the zero.

\item \underline{\textbf{Dropping zeroes:}} Let $S$ be a finite set. Let
$a_{s}$ be an element of $\mathbb{A}$ for each $s\in S$. Let $T$ be a subset
of $S$ such that every $s\in T$ satisfies $a_{s}=0$. Then,%
\begin{equation}
\sum_{s\in S}a_{s}=\sum_{s\in S\setminus T}a_{s}. \label{eq.sum.drop0}%
\end{equation}
(That is, any addends which are zero can be removed from a sum without
changing the sum's value.) See Corollary \ref{cor.ind.gen-com.drop0} below for
a proof of (\ref{eq.sum.drop0}).

\item \underline{\textbf{Renaming the index:}} Let $S$ be a finite set. Let
$a_{s}$ be an element of $\mathbb{A}$ for each $s\in S$. Then,%
\[
\sum_{s\in S}a_{s}=\sum_{t\in S}a_{t}.
\]
This is just saying that the summation index in a sum can be renamed at will,
as long as its name does not clash with other notation.

\item \underline{\textbf{Substituting the index I:}} Let $S$ and $T$ be two
finite sets. Let $f:S\rightarrow T$ be a \textbf{bijective} map. Let $a_{t}$
be an element of $\mathbb{A}$ for each $t\in T$. Then,%
\begin{equation}
\sum_{t\in T}a_{t}=\sum_{s\in S}a_{f\left(  s\right)  }. \label{eq.sum.subs1}%
\end{equation}
(The idea here is that the sum $\sum_{s\in S}a_{f\left(  s\right)  }$ contains
the same addends as the sum $\sum_{t\in T}a_{t}$.) A rigorous proof of
(\ref{eq.sum.subs1}) can be found in Theorem \ref{thm.ind.gen-com.subst1} below.

\textbf{Examples:}

\begin{itemize}
\item For any $n\in\mathbb{N}$, we have%
\[
\sum_{t\in\left\{  1,2,\ldots,n\right\}  }t^{3}=\sum_{s\in\left\{
-n,-n+1,\ldots,-1\right\}  }\left(  -s\right)  ^{3}.
\]
(This follows from (\ref{eq.sum.subs1}), applied to $S=\left\{  -n,-n+1,\ldots
,-1\right\}  $, \newline$T=\left\{  1,2,\ldots,n\right\}  $, $f\left(
s\right)  =-s$, and $a_{t}=t^{3}$.)

\item The sets $S$ and $T$ in (\ref{eq.sum.subs1}) may well be the same. For
example, for any $n\in\mathbb{N}$, we have%
\[
\sum_{t\in\left\{  1,2,\ldots,n\right\}  }t^{3}=\sum_{s\in\left\{
1,2,\ldots,n\right\}  }\left(  n+1-s\right)  ^{3}.
\]
(This follows from (\ref{eq.sum.subs1}), applied to $S=\left\{  1,2,\ldots
,n\right\}  $, $T=\left\{  1,2,\ldots,n\right\}  $, $f\left(  s\right)
=n+1-s$ and $a_{t}=t^{3}$.)

\item More generally: Let $u$ and $v$ be two integers. Then, the map
\newline$\left\{  u,u+1,\ldots,v\right\}  \rightarrow\left\{  u,u+1,\ldots
,v\right\}  $ sending each $s\in\left\{  u,u+1,\ldots,v\right\}  $ to $u+v-s$
is a bijection\footnote{Check this!}. Hence, we can substitute $u+v-s$ for $s$
in the sum $\sum_{s=u}^{v}a_{s}$ whenever an element $a_{s}$ of $\mathbb{A}$
is given for each $s\in\left\{  u,u+1,\ldots,v\right\}  $. We thus obtain the
formula%
\[
\sum_{s=u}^{v}a_{s}=\sum_{s=u}^{v}a_{u+v-s}.
\]

\end{itemize}

\textbf{Remark:}

\begin{itemize}
\item When I use (\ref{eq.sum.subs1}) to rewrite the sum $\sum_{t\in T}a_{t}$
as $\sum_{s\in S}a_{f\left(  s\right)  }$, I say that I have \textquotedblleft
substituted $f\left(  s\right)  $ for $t$ in the sum\textquotedblright.
Conversely, when I use (\ref{eq.sum.subs1}) to rewrite the sum $\sum_{s\in
S}a_{f\left(  s\right)  }$ as $\sum_{t\in T}a_{t}$, I say that I have
\textquotedblleft substituted $t$ for $f\left(  s\right)  $ in the
sum\textquotedblright.

\item For convenience, I have chosen $s$ and $t$ as summation indices in
(\ref{eq.sum.subs1}). But as before, they can be chosen to be any letters not
otherwise used. It is perfectly okay to use one and the same letter for both
of them, e.g., to write $\sum_{s\in T}a_{s}=\sum_{s\in S}a_{f\left(  s\right)
}$.

\item Here is the probably most famous example of substitution in a sum: Fix a
nonnegative integer $n$. Then, we can substitute $n-i$ for $i$ in the sum
$\sum_{i=0}^{n}i$ (since the map $\left\{  0,1,\ldots,n\right\}
\rightarrow\left\{  0,1,\ldots,n\right\}  ,\ i\mapsto n-i$ is a bijection).
Thus, we obtain%
\[
\sum_{i=0}^{n}i=\sum_{i=0}^{n}\left(  n-i\right)  .
\]
Now,%
\begin{align*}
2\sum_{i=0}^{n}i  &  =\sum_{i=0}^{n}i+\underbrace{\sum_{i=0}^{n}i}%
_{=\sum_{i=0}^{n}\left(  n-i\right)  }\ \ \ \ \ \ \ \ \ \ \left(  \text{since
}2q=q+q\text{ for every }q\in\mathbb{Q}\right) \\
&  =\sum_{i=0}^{n}i+\sum_{i=0}^{n}\left(  n-i\right) \\
&  =\sum_{i=0}^{n}\underbrace{\left(  i+\left(  n-i\right)  \right)  }%
_{=n}\ \ \ \ \ \ \ \ \ \ \left(  \text{here, we have used
(\ref{eq.sum.linear1.c}) backwards}\right) \\
&  =\sum_{i=0}^{n}n=\left(  n+1\right)  n\ \ \ \ \ \ \ \ \ \ \left(  \text{by
(\ref{eq.sum.equal})}\right) \\
&  =n\left(  n+1\right)  ,
\end{align*}
and therefore%
\begin{equation}
\sum_{i=0}^{n}i=\dfrac{n\left(  n+1\right)  }{2}. \label{eq.sum.littlegauss1}%
\end{equation}
Since $\sum_{i=0}^{n}i=0+\sum_{i=1}^{n}i=\sum_{i=1}^{n}i$, this rewrites as%
\begin{equation}
\sum_{i=1}^{n}i=\dfrac{n\left(  n+1\right)  }{2}. \label{eq.sum.littlegauss2}%
\end{equation}
This is the famous \textquotedblleft Little Gauss formula\textquotedblright%
\ (supposedly discovered by Carl Friedrich Gauss in primary school, but
already known to the Pythagoreans).
\end{itemize}

\item \underline{\textbf{Substituting the index II:}} Let $S$ and $T$ be two
finite sets. Let $f:S\rightarrow T$ be a \textbf{bijective} map. Let $a_{s}$
be an element of $\mathbb{A}$ for each $s\in S$. Then,%
\begin{equation}
\sum_{s\in S}a_{s}=\sum_{t\in T}a_{f^{-1}\left(  t\right)  }.
\label{eq.sum.subs2}%
\end{equation}
This is, of course, just (\ref{eq.sum.subs1}) but applied to $T$, $S$ and
$f^{-1}$ instead of $S$, $T$ and $f$. (Nevertheless, I prefer to mention
(\ref{eq.sum.subs2}) separately because it often is used in this very form.)

\item \underline{\textbf{Telescoping sums:}} Let $u$ and $v$ be two integers
such that $u-1\leq v$. Let $a_{s}$ be an element of $\mathbb{A}$ for each
$s\in\left\{  u-1,u,\ldots,v\right\}  $. Then,%
\begin{equation}
\sum_{s=u}^{v}\left(  a_{s}-a_{s-1}\right)  =a_{v}-a_{u-1}.
\label{eq.sum.telescope}%
\end{equation}

\textbf{Examples:}

\begin{itemize}
\item Let us give a new proof of (\ref{eq.sum.littlegauss2}). Indeed, fix a
nonnegative integer $n$. An easy computation reveals that%
\begin{equation}
s=\dfrac{s\left(  s+1\right)  }{2}-\dfrac{\left(  s-1\right)  \left(  \left(
s-1\right)  +1\right)  }{2} \label{eq.sum.littlegauss2.pf2.1}%
\end{equation}
for each $s\in\mathbb{Z}$. Thus,%
\begin{align*}
\sum_{i=1}^{n}i  &  =\sum_{s=1}^{n}s=\sum_{s=1}^{n}\left(  \dfrac{s\left(
s+1\right)  }{2}-\dfrac{\left(  s-1\right)  \left(  \left(  s-1\right)
+1\right)  }{2}\right)  \ \ \ \ \ \ \ \ \ \ \left(  \text{by
(\ref{eq.sum.littlegauss2.pf2.1})}\right) \\
&  =\dfrac{n\left(  n+1\right)  }{2}-\underbrace{\dfrac{\left(  1-1\right)
\left(  \left(  1-1\right)  +1\right)  }{2}}_{=0}\\
&  \ \ \ \ \ \ \ \ \ \ \left(  \text{by (\ref{eq.sum.telescope}), applied to
}u=1\text{, }v=n\text{ and }a_{s}=\dfrac{s\left(  s+1\right)  }{2}\right) \\
&  =\dfrac{n\left(  n+1\right)  }{2}.
\end{align*}
Thus, (\ref{eq.sum.littlegauss2}) is proven again. This kind of proof works
often when we need to prove a formula like (\ref{eq.sum.littlegauss2}); the
only tricky part was to \textquotedblleft guess\textquotedblright\ the right
value of $a_{s}$, which is straightforward if you know what you are looking
for (you want $a_{n}-a_{0}$ to be $\dfrac{n\left(  n+1\right)  }{2}$), but
rather tricky if you don't.

\item Another application of (\ref{eq.sum.telescope}) is a proof of the
well-known formula%
\[
\sum_{i=1}^{n}i^{2}=\dfrac{n\left(  n+1\right)  \left(  2n+1\right)  }%
{6}\ \ \ \ \ \ \ \ \ \ \text{for all }n\in\mathbb{N}.
\]
Indeed, an easy computation reveals that%
\[
s^{2}=\dfrac{s\left(  s+1\right)  \left(  2s+1\right)  }{6}-\dfrac{\left(
s-1\right)  \left(  \left(  s-1\right)  +1\right)  \left(  2\left(
s-1\right)  +1\right)  }{6}%
\]
for each $s\in\mathbb{Z}$; now, as in the previous example, we can sum this
equality over all $s\in\left\{  1,2,\ldots,n\right\}  $ and apply
(\ref{eq.sum.telescope}) to obtain our claim.

\item Here is another important identity that follows from
(\ref{eq.sum.telescope}): If $a$ and $b$ are any elements of $\mathbb{A}$, and
if $m\in\mathbb{N}$, then%
\begin{equation}
\left(  a-b\right)  \sum_{i=0}^{m-1}a^{i}b^{m-1-i}=a^{m}-b^{m}.
\label{eq.sum.geometric-sum}%
\end{equation}
(This is one of the versions of the \textquotedblleft geometric series
formula\textquotedblright.) To prove (\ref{eq.sum.geometric-sum}), we observe
that%
\begin{align*}
&  \left(  a-b\right)  \sum_{i=0}^{m-1}a^{i}b^{m-1-i}\\
&  =\sum_{i=0}^{m-1}\underbrace{\left(  a-b\right)  a^{i}b^{m-1-i}}%
_{=aa^{i}b^{m-1-i}-ba^{i}b^{m-1-i}}\ \ \ \ \ \ \ \ \ \ \left(  \text{this
follows from (\ref{eq.sum.linear2})}\right) \\
&  =\sum_{i=0}^{m-1}\left(  aa^{i}b^{m-1-i}-\underbrace{ba^{i}}_{=a^{i}%
b}b^{m-1-i}\right) \\
&  =\sum_{i=0}^{m-1}\left(  \underbrace{aa^{i}}_{=a^{i+1}}b^{m-1-i}%
-\underbrace{a^{i}}_{\substack{=a^{\left(  i-1\right)  +1}\\\text{(since
}i=\left(  i-1\right)  +1\text{)}}}\ \ \underbrace{bb^{m-1-i}}%
_{\substack{=b^{\left(  m-1-i\right)  +1}=b^{m-1-\left(  i-1\right)
}\\\text{(since }\left(  m-1-i\right)  +1=m-1-\left(  i-1\right)  \text{)}%
}}\right) \\
&  =\sum_{i=0}^{m-1}\left(  a^{i+1}b^{m-1-i}-a^{\left(  i-1\right)
+1}b^{m-1-\left(  i-1\right)  }\right) \\
&  =\sum_{s=0}^{m-1}\left(  a^{s+1}b^{m-1-s}-a^{\left(  s-1\right)
+1}b^{m-1-\left(  s-1\right)  }\right) \\
&  \ \ \ \ \ \ \ \ \ \ \left(  \text{here, we have renamed the summation index
}i\text{ as }s\right) \\
&  =\underbrace{a^{\left(  m-1\right)  +1}}_{\substack{=a^{m}\\\text{(since
}\left(  m-1\right)  +1=m\text{)}}}\ \ \underbrace{b^{m-1-\left(  m-1\right)
}}_{\substack{=b^{0}\\\text{(since }m-1-\left(  m-1\right)  =0\text{)}%
}}-\underbrace{a^{\left(  0-1\right)  +1}}_{\substack{=a^{0}\\\text{(since
}\left(  0-1\right)  +1=0\text{)}}}\ \ \underbrace{b^{m-1-\left(  0-1\right)
}}_{\substack{=b^{m}\\\text{(since }m-1-\left(  0-1\right)  =m\text{)}}}\\
&  \ \ \ \ \ \ \ \ \ \ \left(  \text{by (\ref{eq.sum.telescope}) (applied to
}u=0\text{, }v=m-1\text{ and }a_{s}=a^{s+1}b^{m-1-s}\text{)}\right) \\
&  =a^{m}\underbrace{b^{0}}_{=1}-\underbrace{a^{0}}_{=1}b^{m}=a^{m}-b^{m}.
\end{align*}

\item Other examples for the use of (\ref{eq.sum.telescope}) can be found on
\href{https://en.wikipedia.org/wiki/Telescoping_series}{the Wikipedia page for
\textquotedblleft telescoping series\textquotedblright}. Let me add just one
more example: Given $n\in\mathbb{N}$, we want to compute $\sum_{i=1}^{n}%
\dfrac{1}{\sqrt{i}+\sqrt{i+1}}$. (Here, of course, we need to take
$\mathbb{A}=\mathbb{R}$ or $\mathbb{A}=\mathbb{C}$.) We proceed as follows:
For every positive integer $i$, we have%
\[
\dfrac{1}{\sqrt{i}+\sqrt{i+1}}=\dfrac{\left(  \sqrt{i+1}-\sqrt{i}\right)
}{\left(  \sqrt{i}+\sqrt{i+1}\right)  \left(  \sqrt{i+1}-\sqrt{i}\right)
}=\sqrt{i+1}-\sqrt{i}%
\]
(since $\left(  \sqrt{i}+\sqrt{i+1}\right)  \left(  \sqrt{i+1}-\sqrt
{i}\right)  =\left(  \sqrt{i+1}\right)  ^{2}-\left(  \sqrt{i}\right)
^{2}=\left(  i+1\right)  -i=1$). Thus,%
\begin{align*}
&  \sum_{i=1}^{n}\dfrac{1}{\sqrt{i}+\sqrt{i+1}}\\
&  =\sum_{i=1}^{n}\left(  \sqrt{i+1}-\sqrt{i}\right)  =\sum_{s=2}^{n+1}\left(
\sqrt{s}-\sqrt{s-1}\right) \\
&  \ \ \ \ \ \ \ \ \ \ \left(
\begin{array}
[c]{c}%
\text{here, we have substituted }s-1\text{ for }i\text{ in the sum,}\\
\text{since the map }\left\{  2,3,\ldots,n+1\right\}  \rightarrow\left\{
1,2,\ldots,n\right\}  ,\ s\mapsto s-1\\
\text{is a bijection}%
\end{array}
\right) \\
&  =\sqrt{n+1}-\underbrace{\sqrt{2-1}}_{=\sqrt{1}=1}\\
&  \ \ \ \ \ \ \ \ \ \ \left(  \text{by (\ref{eq.sum.telescope}), applied to
}u=2\text{, }v=n+1\text{ and }a_{s}=\sqrt{s}-\sqrt{s-1}\right) \\
&  =\sqrt{n+1}-1.
\end{align*}

\end{itemize}

\textbf{Remarks:}

\begin{itemize}
\item When we use the equality (\ref{eq.sum.telescope}) to rewrite the sum
$\sum_{s=u}^{v}\left(  a_{s}-a_{s-1}\right)  $ as $a_{v}-a_{u-1}$, we can say
that the sum $\sum_{s=u}^{v}\left(  a_{s}-a_{s-1}\right)  $ \textquotedblleft
telescopes\textquotedblright\ to $a_{v}-a_{u-1}$. A sum like $\sum_{s=u}%
^{v}\left(  a_{s}-a_{s-1}\right)  $ is said to be a \textquotedblleft
telescoping sum\textquotedblright. This terminology references the idea that
the sum $\sum_{s=u}^{v}\left(  a_{s}-a_{s-1}\right)  $ \textquotedblleft
shrink\textquotedblright\ to the simple difference $a_{v}-a_{u-1}$ like a
telescope does when it is collapsed.

\item Here is a \textit{proof of (\ref{eq.sum.telescope}):} Let $u$ and $v$ be
two integers such that $u-1\leq v$. Let $a_{s}$ be an element of $\mathbb{A}$
for each $s\in\left\{  u-1,u,\ldots,v\right\}  $. Then,
(\ref{eq.sum.linear1.c}) (applied to $a_{s}-a_{s-1}$ and $a_{s-1}$ instead of
$a_{s}$ and $b_{s}$) yields%
\[
\sum_{s=u}^{v}\left(  \left(  a_{s}-a_{s-1}\right)  +a_{s-1}\right)
=\sum_{s=u}^{v}\left(  a_{s}-a_{s-1}\right)  +\sum_{s=u}^{v}a_{s-1}.
\]
Solving this equation for $\sum_{s=u}^{v}\left(  a_{s}-a_{s-1}\right)  $, we
obtain%
\begin{align}
\sum_{s=u}^{v}\left(  a_{s}-a_{s-1}\right)   &  =\sum_{s=u}^{v}%
\underbrace{\left(  \left(  a_{s}-a_{s-1}\right)  +a_{s-1}\right)  }_{=a_{s}%
}-\underbrace{\sum_{s=u}^{v}a_{s-1}}_{\substack{=\sum_{s=u-1}^{v-1}%
a_{s}\\\text{(here, we have substituted }s\text{ for }s-1\\\text{in the sum)}%
}}\nonumber\\
&  =\sum_{s=u}^{v}a_{s}-\sum_{s=u-1}^{v-1}a_{s}. \label{pf.eq.sum.telescope.1}%
\end{align}

But $u-1\leq v$. Hence, we can split off the addend for $s=u-1$ from the sum
$\sum_{s=u-1}^{v}a_{s}$. We thus obtain
\[
\sum_{s=u-1}^{v}a_{s}=a_{u-1}+\sum_{s=u}^{v}a_{s}.
\]
Solving this equation for $\sum_{s=u}^{v}a_{s}$, we obtain%
\begin{equation}
\sum_{s=u}^{v}a_{s}=\sum_{s=u-1}^{v}a_{s}-a_{u-1}.
\label{pf.eq.sum.telescope.2}%
\end{equation}

Also, $u-1\leq v$. Hence, we can split off the addend for $s=v$ from the sum
$\sum_{s=u-1}^{v}a_{s}$. We thus obtain
\[
\sum_{s=u-1}^{v}a_{s}=a_{v}+\sum_{s=u-1}^{v-1}a_{s}.
\]
Solving this equation for $\sum_{s=u-1}^{v-1}a_{s}$, we obtain%
\begin{equation}
\sum_{s=u-1}^{v-1}a_{s}=\sum_{s=u-1}^{v}a_{s}-a_{v}.
\label{pf.eq.sum.telescope.3}%
\end{equation}
Now, (\ref{pf.eq.sum.telescope.1}) becomes%
\begin{align*}
\sum_{s=u}^{v}\left(  a_{s}-a_{s-1}\right)   &  =\underbrace{\sum_{s=u}%
^{v}a_{s}}_{\substack{=\sum_{s=u-1}^{v}a_{s}-a_{u-1}\\\text{(by
(\ref{pf.eq.sum.telescope.2}))}}}-\underbrace{\sum_{s=u-1}^{v-1}a_{s}%
}_{\substack{=\sum_{s=u-1}^{v}a_{s}-a_{v}\\\text{(by
(\ref{pf.eq.sum.telescope.3}))}}}\\
&  =\left(  \sum_{s=u-1}^{v}a_{s}-a_{u-1}\right)  -\left(  \sum_{s=u-1}%
^{v}a_{s}-a_{v}\right)  =a_{v}-a_{u-1}.
\end{align*}
This proves (\ref{eq.sum.telescope}).
\end{itemize}

\item \underline{\textbf{Restricting to a subset:}} Let $S$ be a finite set.
Let $T$ be a subset of $S$. Let $a_{s}$ be an element of $\mathbb{A}$ for each
$s\in T$. Then,%
\[
\sum_{\substack{s\in S;\\s\in T}}a_{s}=\sum_{s\in T}a_{s}.
\]
This is because the $s\in S$ satisfying $s\in T$ are exactly the elements of
$T$.

\textbf{Remark:} Here is a slightly more general form of this rule: Let $S$ be
a finite set. Let $T$ be a subset of $S$. Let $\mathcal{A}\left(  s\right)  $
be a logical statement for each $s\in S$. Let $a_{s}$ be an element of
$\mathbb{A}$ for each $s\in T$ satisfying $\mathcal{A}\left(  s\right)  $.
Then,%
\[
\sum_{\substack{s\in S;\\s\in T;\\\mathcal{A}\left(  s\right)  }}a_{s}%
=\sum_{\substack{s\in T;\\\mathcal{A}\left(  s\right)  }}a_{s}.
\]

\item \underline{\textbf{Splitting a sum by a value of a function:}} Let $S$
be a finite set. Let $W$ be a set. Let $f:S\rightarrow W$ be a map. Let
$a_{s}$ be an element of $\mathbb{A}$ for each $s\in S$. Then,%
\begin{equation}
\sum_{s\in S}a_{s}=\sum_{w\in W}\sum_{\substack{s\in S;\\f\left(  s\right)
=w}}a_{s}. \label{eq.sum.sheph}%
\end{equation}
The idea behind this formula is the following: The left hand side is the sum
of all $a_{s}$ for $s\in S$. The right hand side is the same sum, but split in
a particular way: First, for each $w\in W$, we sum the $a_{s}$ for all $s\in
S$ satisfying $f\left(  s\right)  =w$, and then we take the sum of all these
\textquotedblleft partial sums\textquotedblright. For a rigorous proof of
(\ref{eq.sum.sheph}), see Theorem \ref{thm.ind.gen-com.shephf} (for the case
when $W$ is finite) and Theorem \ref{thm.ind.gen-com.sheph} (for the general case).

\textbf{Examples:}

\begin{itemize}
\item Let $n\in\mathbb{N}$. Then,%
\begin{equation}
\sum_{s\in\left\{  -n,-\left(  n-1\right)  ,\ldots,n\right\}  }s^{3}%
=\sum_{w\in\left\{  0,1,\ldots,n\right\}  }\sum_{\substack{s\in\left\{
-n,-\left(  n-1\right)  ,\ldots,n\right\}  ;\\\left\vert s\right\vert
=w}}s^{3}. \label{eq.sum.sheph.exam1}%
\end{equation}
(This follows from (\ref{eq.sum.sheph}), applied to $S=\left\{  -n,-\left(
n-1\right)  ,\ldots,n\right\}  $, $W=\left\{  0,1,\ldots,n\right\}  $ and
$f\left(  s\right)  =\left\vert s\right\vert $.) You might wonder what you
gain by this observation. But actually, it allows you to compute the sum: For
any $w\in\left\{  0,1,\ldots,n\right\}  $, the sum $\sum_{\substack{s\in
\left\{  -n,-\left(  n-1\right)  ,\ldots,n\right\}  ;\\\left\vert s\right\vert
=w}}s^{3}$ is $0$\ \ \ \ \footnote{\textit{Proof.} If $w=0$, then this sum
$\sum_{\substack{s\in\left\{  -n,-\left(  n-1\right)  ,\ldots,n\right\}
;\\\left\vert s\right\vert =w}}s^{3}$ consists of one addend only, and this
addend is $0^{3}$. If $w>0$, then this sum has two addends, namely $\left(
-w\right)  ^{3}$ and $w^{3}$. In either case, the sum is $0$ (because
$0^{3}=0$ and $\left(  -w\right)  ^{3}+w^{3}=-w^{3}+w^{3}=0$).}, and therefore
(\ref{eq.sum.sheph.exam1}) becomes%
\[
\sum_{s\in\left\{  -n,-\left(  n-1\right)  ,\ldots,n\right\}  }s^{3}%
=\sum_{w\in\left\{  0,1,\ldots,n\right\}  }\ \ \underbrace{\sum
_{\substack{s\in\left\{  -n,-\left(  n-1\right)  ,\ldots,n\right\}
;\\\left\vert s\right\vert =w}}s^{3}}_{=0}=\sum_{w\in\left\{  0,1,\ldots
,n\right\}  }0=0.
\]
Thus, a strategic application of (\ref{eq.sum.sheph}) can help in evaluating a sum.

\item Let $S$ be a finite set. Let $W$ be a set. Let $f:S\rightarrow W$ be a
map. If we apply (\ref{eq.sum.sheph}) to $a_{s}=1$, then we obtain%
\[
\sum_{s\in S}1=\sum_{w\in W}\underbrace{\sum_{\substack{s\in S;\\f\left(
s\right)  =w}}1}_{\substack{=\left\vert \left\{  s\in S\ \mid\ f\left(
s\right)  =w\right\}  \right\vert \cdot1\\=\left\vert \left\{  s\in
S\ \mid\ f\left(  s\right)  =w\right\}  \right\vert }}=\sum_{w\in W}\left\vert
\left\{  s\in S\ \mid\ f\left(  s\right)  =w\right\}  \right\vert .
\]
Since $\sum_{s\in S}1=\left\vert S\right\vert \cdot1=\left\vert S\right\vert
$, this rewrites as follows:%
\begin{equation}
\left\vert S\right\vert =\sum_{w\in W}\left\vert \left\{  s\in S\ \mid
\ f\left(  s\right)  =w\right\}  \right\vert . \label{eq.sum.sheph.exam2}%
\end{equation}
This equality is often called the \textit{shepherd's principle}, because it is
connected to the joke that \textquotedblleft in order to count a flock of
sheep, just count the legs and divide by $4$\textquotedblright. The connection
is somewhat weak, actually; the equality (\ref{eq.sum.sheph.exam2}) is better
regarded as a formalization of the (less funny) idea that in order to count
all legs of a flock of sheep, you can count the legs of every single sheep,
and then sum the resulting numbers over all sheep in the flock. Think of the
$S$ in (\ref{eq.sum.sheph.exam2}) as the set of all legs of all sheep in the
flock; think of $W$ as the set of all sheep in the flock; and think of $f$ as
the function which sends every leg to the (hopefully uniquely determined)
sheep it belongs to.
\end{itemize}

\textbf{Remarks:}

\begin{itemize}
\item If $f:S\rightarrow W$ is a map between two sets $S$ and $W$, and if $w$
is an element of $W$, then it is common to denote the set $\left\{  s\in
S\ \mid\ f\left(  s\right)  =w\right\}  $ by $f^{-1}\left(  w\right)  $.
(Formally speaking, this notation might clash with the notation $f^{-1}\left(
w\right)  $ for the actual preimage of $w$ when $f$ happens to be bijective;
but in practice, this causes far less confusion than it might seem to.) Using
this notation, we can rewrite (\ref{eq.sum.sheph}) as follows:%
\begin{equation}
\sum_{s\in S}a_{s}=\sum_{w\in W}\ \ \underbrace{\sum_{\substack{s\in
S;\\f\left(  s\right)  =w}}}_{=\sum_{s\in f^{-1}\left(  w\right)  }}a_{s}%
=\sum_{w\in W}\ \ \sum_{s\in f^{-1}\left(  w\right)  }a_{s}.
\label{eq.sum.sheph.preimg}%
\end{equation}

\item When I rewrite a sum $\sum_{s\in S}a_{s}$ as $\sum_{w\in W}%
\ \ \sum_{\substack{s\in S;\\f\left(  s\right)  =w}}a_{s}$ (or as $\sum_{w\in
W}\ \ \sum_{s\in f^{-1}\left(  w\right)  }a_{s}$), I say that I am
\textquotedblleft splitting the sum according to the value of $f\left(
s\right)  $\textquotedblright. (Though, most of the time, I shall be doing
such manipulations without explicit mention.)
\end{itemize}

\item \underline{\textbf{Splitting a sum into subsums:}} Let $S$ be a finite
set. Let $S_{1},S_{2},\ldots,S_{n}$ be finitely many subsets of $S$. Assume
that these subsets $S_{1},S_{2},\ldots,S_{n}$ are pairwise disjoint (i.e., we
have $S_{i}\cap S_{j}=\varnothing$ for any two distinct elements $i$ and $j$
of $\left\{  1,2,\ldots,n\right\}  $) and their union is $S$. (Thus, every
element of $S$ lies in precisely one of the subsets $S_{1},S_{2},\ldots,S_{n}%
$.) Let $a_{s}$ be an element of $\mathbb{A}$ for each $s\in S$. Then,%
\begin{equation}
\sum_{s\in S}a_{s}=\sum_{w=1}^{n}\ \ \sum_{s\in S_{w}}a_{s}.
\label{eq.sum.split-n}%
\end{equation}
This is a generalization of (\ref{eq.sum.split}) (indeed, (\ref{eq.sum.split})
is obtained from (\ref{eq.sum.split-n}) by setting $n=2$, $S_{1}=X$ and
$S_{2}=Y$). It is also a consequence of (\ref{eq.sum.sheph}): Indeed, set
$W=\left\{  1,2,\ldots,n\right\}  $, and define a map $f:S\rightarrow W$ to
send each $s\in S$ to the unique $w\in\left\{  1,2,\ldots,n\right\}  $ for
which $s\in S_{w}$. Then, every $w\in W$ satisfies $\sum_{\substack{s\in
S;\\f\left(  s\right)  =w}}a_{s}=\sum_{s\in S_{w}}a_{s}$; therefore,
(\ref{eq.sum.sheph}) becomes (\ref{eq.sum.split-n}).

\textbf{Example:} If we set $a_{s}=1$ for each $s\in S$, then
(\ref{eq.sum.split-n}) becomes%
\[
\sum_{s\in S}1=\sum_{w=1}^{n}\ \ \underbrace{\sum_{s\in S_{w}}1}_{=\left\vert
S_{w}\right\vert }=\sum_{w=1}^{n}\left\vert S_{w}\right\vert .
\]
Hence,
\[
\sum_{w=1}^{n}\left\vert S_{w}\right\vert =\sum_{s\in S}1=\left\vert
S\right\vert \cdot1=\left\vert S\right\vert .
\]

\item \underline{\textbf{Fubini's theorem (interchanging the order of
summation):}} Let $X$ and $Y$ be two finite sets. Let $a_{\left(  x,y\right)
}$ be an element of $\mathbb{A}$ for each $\left(  x,y\right)  \in X\times Y$.
Then,%
\begin{equation}
\sum_{x\in X}\ \ \sum_{y\in Y}a_{\left(  x,y\right)  }=\sum_{\left(
x,y\right)  \in X\times Y}a_{\left(  x,y\right)  }=\sum_{y\in Y}\ \ \sum_{x\in
X}a_{\left(  x,y\right)  }. \label{eq.sum.fubini}%
\end{equation}
This is called \textit{Fubini's theorem for finite sums}, and is a lot easier
to prove than what analysts tend to call Fubini's theorem. I shall sketch a
proof shortly (in the Remark below); but first, let me give some intuition for
the statement. Imagine that you have a rectangular table filled with numbers.
If you want to sum the numbers in the table, you can proceed in several ways.
One way is to sum the numbers in each row, and then sum all the sums you have
obtained. Another way is to sum the numbers in each column, and then sum all
the obtained sums. Either way, you get the same result -- namely, the sum of
all numbers in the table. This is essentially what (\ref{eq.sum.fubini}) says,
at least when $X=\left\{  1,2,\ldots,n\right\}  $ and $Y=\left\{
1,2,\ldots,m\right\}  $ for some integers $n$ and $m$. In this case, the
numbers $a_{\left(  x,y\right)  }$ can be viewed as forming a table, where
$a_{\left(  x,y\right)  }$ is placed in the cell at the intersection of row
$x$ with column $y$. When $X$ and $Y$ are arbitrary finite sets (not
necessarily $\left\{  1,2,\ldots,n\right\}  $ and $\left\{  1,2,\ldots
,m\right\}  $), then you need to slightly stretch your imagination in order to
see the $a_{\left(  x,y\right)  }$ as \textquotedblleft forming a
table\textquotedblright; in fact, there is no obvious order in which the
numbers appear in a row or column, but there is still a notion of rows and columns.

\textbf{Examples:}

\begin{itemize}
\item Let $n\in\mathbb{N}$ and $m\in\mathbb{N}$. Let $a_{\left(  x,y\right)
}$ be an element of $\mathbb{A}$ for each $\left(  x,y\right)  \in\left\{
1,2,\ldots,n\right\}  \times\left\{  1,2,\ldots,m\right\}  $. Then,
\begin{equation}
\sum_{x=1}^{n}\ \ \sum_{y=1}^{m}a_{\left(  x,y\right)  }=\sum_{\left(
x,y\right)  \in\left\{  1,2,\ldots,n\right\}  \times\left\{  1,2,\ldots
,m\right\}  }a_{\left(  x,y\right)  }=\sum_{y=1}^{m}\ \ \sum_{x=1}%
^{n}a_{\left(  x,y\right)  }. \label{eq.sum.fubini.nm}%
\end{equation}
(This follows from (\ref{eq.sum.fubini}), applied to $X=\left\{
1,2,\ldots,n\right\}  $ and $Y=\left\{  1,2,\ldots,m\right\}  $.) We can
rewrite the equality (\ref{eq.sum.fubini.nm}) without using $\sum$ signs; it
then takes the following form:%
\begin{align*}
&  \left(  a_{\left(  1,1\right)  }+a_{\left(  1,2\right)  }+\cdots+a_{\left(
1,m\right)  }\right) \\
&  \ \ \ \ \ \ \ \ \ \ +\left(  a_{\left(  2,1\right)  }+a_{\left(
2,2\right)  }+\cdots+a_{\left(  2,m\right)  }\right) \\
&  \ \ \ \ \ \ \ \ \ \ +\cdots\\
&  \ \ \ \ \ \ \ \ \ \ +\left(  a_{\left(  n,1\right)  }+a_{\left(
n,2\right)  }+\cdots+a_{\left(  n,m\right)  }\right) \\
&  =a_{\left(  1,1\right)  }+a_{\left(  1,2\right)  }+\cdots+a_{\left(
n,m\right)  }\ \ \ \ \ \ \ \ \ \ \left(  \text{this is the sum of all
}nm\text{ numbers }a_{\left(  x,y\right)  }\right) \\
&  =\left(  a_{\left(  1,1\right)  }+a_{\left(  2,1\right)  }+\cdots
+a_{\left(  n,1\right)  }\right) \\
&  \ \ \ \ \ \ \ \ \ \ +\left(  a_{\left(  1,2\right)  }+a_{\left(
2,2\right)  }+\cdots+a_{\left(  n,2\right)  }\right) \\
&  \ \ \ \ \ \ \ \ \ \ +\cdots\\
&  \ \ \ \ \ \ \ \ \ \ +\left(  a_{\left(  1,m\right)  }+a_{\left(
2,m\right)  }+\cdots+a_{\left(  n,m\right)  }\right)  .
\end{align*}
In other words, we can sum the entries of the rectangular table%
\[%
\begin{tabular}
[c]{|cccc|}\hline
$a_{\left(  1,1\right)  }$ & $a_{\left(  1,2\right)  }$ & $\cdots$ &
$a_{\left(  1,m\right)  }$\\
$a_{\left(  2,1\right)  }$ & $a_{\left(  2,2\right)  }$ & $\cdots$ &
$a_{\left(  2,m\right)  }$\\
$\vdots$ & $\vdots$ & $\ddots$ & $\vdots$\\
$a_{\left(  n,1\right)  }$ & $a_{\left(  n,2\right)  }$ & $\cdots$ &
$a_{\left(  n,m\right)  }$\\\hline
\end{tabular}
\
\]
in three different ways:

\begin{itemize}
\item[\textbf{(a)}] row by row (i.e., first summing the entries in each row,
then summing up the $n$ resulting tallies);

\item[\textbf{(b)}] arbitrarily (i.e., just summing all entries of the table
in some arbitrary order);

\item[\textbf{(c)}] column by column (i.e., first summing the entries in each
column, then summing up the $m$ resulting tallies);
\end{itemize}

and each time, we get the same result.

\item Here is a concrete application of (\ref{eq.sum.fubini.nm}): Let
$n\in\mathbb{N}$ and $m\in\mathbb{N}$. We want to compute $\sum_{\left(
x,y\right)  \in\left\{  1,2,\ldots,n\right\}  \times\left\{  1,2,\ldots
,m\right\}  }xy$. (This is the sum of all entries of the $n\times m$
multiplication table.) Applying (\ref{eq.sum.fubini.nm}) to $a_{\left(
x,y\right)  }=xy$, we obtain%
\[
\sum_{x=1}^{n}\ \ \sum_{y=1}^{m}xy=\sum_{\left(  x,y\right)  \in\left\{
1,2,\ldots,n\right\}  \times\left\{  1,2,\ldots,m\right\}  }xy=\sum_{y=1}%
^{m}\ \ \sum_{x=1}^{n}xy.
\]
Hence,%
\begin{align*}
&  \sum_{\left(  x,y\right)  \in\left\{  1,2,\ldots,n\right\}  \times\left\{
1,2,\ldots,m\right\}  }xy\\
&  =\sum_{x=1}^{n}\ \ \underbrace{\sum_{y=1}^{m}xy}_{\substack{=\sum_{s=1}%
^{m}xs=x\sum_{s=1}^{m}s\\\text{(by (\ref{eq.sum.linear2}), applied to
}S=\left\{  1,2,\ldots,m\right\}  \text{,}\\a_{s}=s\text{ and }\lambda
=x\text{)}}}=\sum_{x=1}^{n}x\underbrace{\sum_{s=1}^{m}s}_{\substack{=\sum
_{i=1}^{m}i=\dfrac{m\left(  m+1\right)  }{2}\\\text{(by
(\ref{eq.sum.littlegauss2}), applied to }m\\\text{instead of }n\text{)}}}\\
&  =\sum_{x=1}^{n}x\dfrac{m\left(  m+1\right)  }{2}=\sum_{x=1}^{n}%
\dfrac{m\left(  m+1\right)  }{2}x=\sum_{s=1}^{n}\dfrac{m\left(  m+1\right)
}{2}s\\
&  =\dfrac{m\left(  m+1\right)  }{2}\underbrace{\sum_{s=1}^{n}s}%
_{\substack{=\sum_{i=1}^{n}i=\dfrac{n\left(  n+1\right)  }{2}\\\text{(by
(\ref{eq.sum.littlegauss2}))}}}\\
&  \ \ \ \ \ \ \ \ \ \ \left(  \text{by (\ref{eq.sum.linear2}), applied to
}S=\left\{  1,2,\ldots,n\right\}  \text{, }a_{s}=s\text{ and }\lambda
=\dfrac{m\left(  m+1\right)  }{2}\right) \\
&  =\dfrac{m\left(  m+1\right)  }{2}\cdot\dfrac{n\left(  n+1\right)  }{2}.
\end{align*}

\end{itemize}

\textbf{Remarks:}

\begin{itemize}
\item I have promised to outline a proof of (\ref{eq.sum.fubini}). Here it
comes: Let $S=X\times Y$ and $W=Y$, and let $f:S\rightarrow W$ be the map
which sends every pair $\left(  x,y\right)  $ to its second entry $y$. Then,
(\ref{eq.sum.sheph.preimg}) shows that%
\begin{equation}
\sum_{s\in X\times Y}a_{s}=\sum_{w\in Y}\sum_{s\in f^{-1}\left(  w\right)
}a_{s}. \label{eq.sum.fubini.pf.1}%
\end{equation}
However, for every given $w\in Y$, the set $f^{-1}\left(  w\right)  $ is
simply the set of all pairs $\left(  x,w\right)  $ with $x\in X$. Thus, for
every given $w\in Y$, there is a bijection $g_{w}:X\rightarrow f^{-1}\left(
w\right)  $ given by%
\[
g_{w}\left(  x\right)  =\left(  x,w\right)  \ \ \ \ \ \ \ \ \ \ \text{for all
}x\in X.
\]
Hence, for every given $w\in Y$, we can substitute $g_{w}\left(  x\right)  $
for $s$ in the sum $\sum_{s\in f^{-1}\left(  w\right)  }a_{s}$, and thus
obtain%
\[
\sum_{s\in f^{-1}\left(  w\right)  }a_{s}=\sum_{x\in X}\underbrace{a_{g_{w}%
\left(  x\right)  }}_{\substack{=a_{\left(  x,w\right)  }\\\text{(since }%
g_{w}\left(  x\right)  =\left(  x,w\right)  \text{)}}}=\sum_{x\in X}a_{\left(
x,w\right)  }.
\]
Hence, (\ref{eq.sum.fubini.pf.1}) becomes%
\[
\sum_{s\in X\times Y}a_{s}=\sum_{w\in Y}\ \ \underbrace{\sum_{s\in
f^{-1}\left(  w\right)  }a_{s}}_{=\sum_{x\in X}a_{\left(  x,w\right)  }}%
=\sum_{w\in Y}\ \ \sum_{x\in X}a_{\left(  x,w\right)  }=\sum_{y\in Y}%
\ \ \sum_{x\in X}a_{\left(  x,y\right)  }%
\]
(here, we have renamed the summation index $w$ as $y$ in the outer sum).
Therefore,%
\[
\sum_{y\in Y}\ \ \sum_{x\in X}a_{\left(  x,y\right)  }=\sum_{s\in X\times
Y}a_{s}=\sum_{\left(  x,y\right)  \in X\times Y}a_{\left(  x,y\right)  }%
\]
(here, we have renamed the summation index $s$ as $\left(  x,y\right)  $).
Thus, we have proven the second part of the equality (\ref{eq.sum.fubini}).
The first part can be proven similarly.

\item I like to abbreviate the equality (\ref{eq.sum.fubini.nm}) as follows:%
\begin{equation}
\sum_{x=1}^{n}\ \ \sum_{y=1}^{m}=\sum_{\left(  x,y\right)  \in\left\{
1,2,\ldots,n\right\}  \times\left\{  1,2,\ldots,m\right\}  }=\sum_{y=1}%
^{m}\ \ \sum_{x=1}^{n}. \label{eq.sum.fubini.eq-sums}%
\end{equation}
This is an \textquotedblleft equality between summation
signs\textquotedblright; it should be understood as follows: Every time you
see an \textquotedblleft$\sum_{x=1}^{n}\sum_{y=1}^{m}$\textquotedblright\ in
an expression, you can replace it by a \textquotedblleft$\sum_{\left(
x,y\right)  \in\left\{  1,2,\ldots,n\right\}  \times\left\{  1,2,\ldots
,m\right\}  }$\textquotedblright\ or by a \textquotedblleft$\sum_{y=1}^{m}%
\sum_{x=1}^{n}$\textquotedblright, and similarly the other ways round.
\end{itemize}

\item \underline{\textbf{Triangular Fubini's theorem I:}} The equality
(\ref{eq.sum.fubini.nm}) formalizes the idea that we can sum the entries of a
rectangular table by first tallying each row and then adding together, or
first tallying each column and adding together. The same holds for triangular
tables. More precisely: Let $n\in\mathbb{N}$. Let $T_{n}$ be the set $\left\{
\left(  x,y\right)  \in\left\{  1,2,3,\ldots\right\}  ^{2}\ \mid\ x+y\leq
n\right\}  $. (For instance, if $n=3$, then $T_{n}=T_{3}=\left\{  \left(
1,1\right)  ,\left(  1,2\right)  ,\left(  2,1\right)  \right\}  $.) Let
$a_{\left(  x,y\right)  }$ be an element of $\mathbb{A}$ for each $\left(
x,y\right)  \in T_{n}$. Then,%
\begin{equation}
\sum_{x=1}^{n}\ \ \sum_{y=1}^{n-x}a_{\left(  x,y\right)  }=\sum_{\left(
x,y\right)  \in T_{n}}a_{\left(  x,y\right)  }=\sum_{y=1}^{n}\ \ \sum
_{x=1}^{n-y}a_{\left(  x,y\right)  }. \label{eq.sum.fubini.triangle}%
\end{equation}

\textbf{Examples:}

\begin{itemize}
\item In the case when $n=4$, the formula (\ref{eq.sum.fubini.triangle})
(rewritten without the use of $\sum$ signs) looks as follows:%
\begin{align*}
&  \left(  a_{\left(  1,1\right)  }+a_{\left(  1,2\right)  }+a_{\left(
1,3\right)  }\right)  +\left(  a_{\left(  2,1\right)  }+a_{\left(  2,2\right)
}\right)  +a_{\left(  3,1\right)  }+\left(  \text{empty sum}\right) \\
&  =\left(  \text{the sum of the }a_{\left(  x,y\right)  }\text{ for all
}\left(  x,y\right)  \in T_{4}\right) \\
&  =\left(  a_{\left(  1,1\right)  }+a_{\left(  2,1\right)  }+a_{\left(
3,1\right)  }\right)  +\left(  a_{\left(  1,2\right)  }+a_{\left(  2,2\right)
}\right)  +a_{\left(  1,3\right)  }+\left(  \text{empty sum}\right)  .
\end{align*}
In other words, we can sum the entries of the triangular table%
\[%
\begin{tabular}
[c]{|ccc}\hline
$a_{\left(  1,1\right)  }$ & $a_{\left(  1,2\right)  }$ &
\multicolumn{1}{c|}{$a_{\left(  1,3\right)  }$}\\\cline{3-3}%
$a_{\left(  2,1\right)  }$ & $a_{\left(  2,2\right)  }$ &
\multicolumn{1}{|c}{}\\\cline{2-2}%
$a_{\left(  3,1\right)  }$ & \multicolumn{1}{|c}{} & \\\cline{1-1}%
\end{tabular}
\]
in three different ways:

\begin{itemize}
\item[\textbf{(a)}] row by row (i.e., first summing the entries in each row,
then summing up the resulting tallies);

\item[\textbf{(b)}] arbitrarily (i.e., just summing all entries of the table
in some arbitrary order);

\item[\textbf{(c)}] column by column (i.e., first summing the entries in each
column, then summing up the resulting tallies);
\end{itemize}

and each time, we get the same result.

\item Let us use (\ref{eq.sum.fubini.triangle}) to compute $\left\vert
T_{n}\right\vert $. Indeed, we can apply (\ref{eq.sum.fubini.triangle}) to
$a_{\left(  x,y\right)  }=1$. Thus, we obtain%
\[
\sum_{x=1}^{n}\ \ \sum_{y=1}^{n-x}1=\sum_{\left(  x,y\right)  \in T_{n}}%
1=\sum_{y=1}^{n}\ \ \sum_{x=1}^{n-y}1.
\]
Hence,%
\[
\sum_{x=1}^{n}\ \ \sum_{y=1}^{n-x}1=\sum_{\left(  x,y\right)  \in T_{n}%
}1=\left\vert T_{n}\right\vert ,
\]
so that%
\begin{align*}
\left\vert T_{n}\right\vert  &  =\sum_{x=1}^{n}\ \ \underbrace{\sum
_{y=1}^{n-x}1}_{=n-x}=\sum_{x=1}^{n}\left(  n-x\right)  =\sum_{i=0}^{n-1}i\\
&  \ \ \ \ \ \ \ \ \ \ \left(
\begin{array}
[c]{c}%
\text{here, we have substituted }i\text{ for }n-x\text{ in the sum,}\\
\text{since the map }\left\{  1,2,\ldots,n\right\}  \rightarrow\left\{
0,1,\ldots,n-1\right\}  ,\ x\mapsto n-x\\
\text{is a bijection}%
\end{array}
\right) \\
&  =\dfrac{\left(  n-1\right)  \left(  \left(  n-1\right)  +1\right)  }%
{2}\ \ \ \ \ \ \ \ \ \ \left(  \text{by (\ref{eq.sum.littlegauss1}), applied
to }n-1\text{ instead of }n\right) \\
&  =\dfrac{\left(  n-1\right)  n}{2}.
\end{align*}

\end{itemize}

\textbf{Remarks:}

\begin{itemize}
\item The sum $\sum_{\left(  x,y\right)  \in T_{n}}a_{\left(  x,y\right)  }$
in (\ref{eq.sum.fubini.triangle}) can also be rewritten as $\sum
_{\substack{\left(  x,y\right)  \in\left\{  1,2,3,\ldots\right\}
^{2};\\x+y\leq n}}a_{\left(  x,y\right)  }$.

\item Let us prove (\ref{eq.sum.fubini.triangle}). Indeed, the proof will be
very similar to our proof of (\ref{eq.sum.fubini}) above. Let $S=T_{n}$ and
$W=\left\{  1,2,\ldots,n\right\}  $, and let $f:S\rightarrow W$ be the map
which sends every pair $\left(  x,y\right)  $ to its second entry $y$. Then,
(\ref{eq.sum.sheph.preimg}) shows that%
\begin{equation}
\sum_{s\in T_{n}}a_{s}=\sum_{w\in W}\ \ \sum_{s\in f^{-1}\left(  w\right)
}a_{s}. \label{eq.sum.fubini.triangle.pf.1}%
\end{equation}
However, for every given $w\in W$, the set $f^{-1}\left(  w\right)  $ is
simply the set of all pairs $\left(  x,w\right)  $ with $x\in\left\{
1,2,\ldots,n-w\right\}  $. Thus, for every given $w\in W$, there is a
bijection $g_{w}:\left\{  1,2,\ldots,n-w\right\}  \rightarrow f^{-1}\left(
w\right)  $ given by%
\[
g_{w}\left(  x\right)  =\left(  x,w\right)  \ \ \ \ \ \ \ \ \ \ \text{for all
}x\in\left\{  1,2,\ldots,n-w\right\}  .
\]
Hence, for every given $w\in W$, we can substitute $g_{w}\left(  x\right)  $
for $s$ in the sum $\sum_{s\in f^{-1}\left(  w\right)  }a_{s}$, and thus
obtain%
\[
\sum_{s\in f^{-1}\left(  w\right)  }a_{s}=\underbrace{\sum_{x\in\left\{
1,2,\ldots,n-w\right\}  }}_{=\sum_{x=1}^{n-w}}\underbrace{a_{g_{w}\left(
x\right)  }}_{\substack{=a_{\left(  x,w\right)  }\\\text{(since }g_{w}\left(
x\right)  =\left(  x,w\right)  \text{)}}}=\sum_{x=1}^{n-w}a_{\left(
x,w\right)  }.
\]
Hence, (\ref{eq.sum.fubini.triangle.pf.1}) becomes%
\[
\sum_{s\in T_{n}}a_{s}=\underbrace{\sum_{w\in W}}_{\substack{=\sum_{w=1}%
^{n}\\\text{(since }W=\left\{  1,2,\ldots,n\right\}  \text{)}}%
}\underbrace{\sum_{s\in f^{-1}\left(  w\right)  }a_{s}}_{=\sum_{x=1}%
^{n-w}a_{\left(  x,w\right)  }}=\sum_{w=1}^{n}\ \ \sum_{x=1}^{n-w}a_{\left(
x,w\right)  }=\sum_{y=1}^{n}\ \ \sum_{x=1}^{n-y}a_{\left(  x,y\right)  }%
\]
(here, we have renamed the summation index $w$ as $y$ in the outer sum).
Therefore,%
\[
\sum_{y=1}^{n}\ \ \sum_{x=1}^{n-y}a_{\left(  x,y\right)  }=\sum_{s\in T_{n}%
}a_{s}=\sum_{\left(  x,y\right)  \in T_{n}}a_{\left(  x,y\right)  }.
\]
Thus, we have proven the second part of the equality
(\ref{eq.sum.fubini.triangle}). The first part can be proven similarly.
\end{itemize}

\item \underline{\textbf{Triangular Fubini's theorem II:}} Here is another
equality similar to (\ref{eq.sum.fubini.triangle}). Let $n\in\mathbb{N}$. Let
$Q_{n}$ be the set $\left\{  \left(  x,y\right)  \in\left\{  1,2,\ldots
,n\right\}  ^{2}\ \mid\ x\leq y\right\}  $. (For instance, if $n=3$, then
$Q_{n}=Q_{3}=\left\{  \left(  1,1\right)  ,\left(  1,2\right)  ,\left(
1,3\right)  ,\left(  2,2\right)  ,\left(  2,3\right)  ,\left(  3,3\right)
\right\}  $.) Let $a_{\left(  x,y\right)  }$ be an element of $\mathbb{A}$ for
each $\left(  x,y\right)  \in Q_{n}$. Then,%
\begin{equation}
\sum_{x=1}^{n}\ \ \sum_{y=x}^{n}a_{\left(  x,y\right)  }=\sum_{\left(
x,y\right)  \in Q_{n}}a_{\left(  x,y\right)  }=\sum_{y=1}^{n}\ \ \sum
_{x=1}^{y}a_{\left(  x,y\right)  }. \label{eq.sum.fubini.triangle2}%
\end{equation}

\textbf{Examples:}

\begin{itemize}
\item In the case when $n=4$, the formula (\ref{eq.sum.fubini.triangle2})
(rewritten without the use of $\sum$ signs) looks as follows:%
\begin{align*}
&  \left(  a_{\left(  1,1\right)  }+a_{\left(  1,2\right)  }+a_{\left(
1,3\right)  }+a_{\left(  1,4\right)  }\right) \\
&  \ \ \ \ \ \ \ \ \ \ +\left(  a_{\left(  2,2\right)  }+a_{\left(
2,3\right)  }+a_{\left(  2,4\right)  }\right) \\
&  \ \ \ \ \ \ \ \ \ \ +\left(  a_{\left(  3,3\right)  }+a_{\left(
3,4\right)  }\right) \\
&  \ \ \ \ \ \ \ \ \ \ +a_{\left(  4,4\right)  }\\
&  =\left(  \text{the sum of the }a_{\left(  x,y\right)  }\text{ for all
}\left(  x,y\right)  \in Q_{4}\right) \\
&  =a_{\left(  1,1\right)  }\\
&  \ \ \ \ \ \ \ \ \ \ +\left(  a_{\left(  1,2\right)  }+a_{\left(
2,2\right)  }\right) \\
&  \ \ \ \ \ \ \ \ \ \ +\left(  a_{\left(  1,3\right)  }+a_{\left(
2,3\right)  }+a_{\left(  3,3\right)  }\right) \\
&  \ \ \ \ \ \ \ \ \ \ +\left(  a_{\left(  1,4\right)  }+a_{\left(
2,4\right)  }+a_{\left(  3,4\right)  }+a_{\left(  4,4\right)  }\right)  .
\end{align*}
In other words, we can sum the entries of the triangular table%
\[%
\begin{tabular}
[c]{cccc|}\hline
\multicolumn{1}{|c}{$a_{\left(  1,1\right)  }$} & $a_{\left(  1,2\right)  }$ &
$a_{\left(  1,3\right)  }$ & $a_{\left(  1,4\right)  }$\\\cline{1-1}
& \multicolumn{1}{|c}{$a_{\left(  2,2\right)  }$} & $a_{\left(  2,3\right)  }$
& $a_{\left(  2,4\right)  }$\\\cline{2-2}
&  & \multicolumn{1}{|c}{$a_{\left(  3,3\right)  }$} & $a_{\left(  3,4\right)
}$\\\cline{3-3}
&  &  & \multicolumn{1}{|c|}{$a_{\left(  4,4\right)  }$}\\\cline{4-4}%
\end{tabular}
\]
in three different ways:

\begin{itemize}
\item[\textbf{(a)}] row by row (i.e., first summing the entries in each row,
then summing up the resulting tallies);

\item[\textbf{(b)}] arbitrarily (i.e., just summing all entries of the table
in some arbitrary order);

\item[\textbf{(c)}] column by column (i.e., first summing the entries in each
column, then summing up the resulting tallies);
\end{itemize}

and each time, we get the same result.

\item Let us use (\ref{eq.sum.fubini.triangle2}) to compute $\left\vert
Q_{n}\right\vert $. Indeed, we can apply (\ref{eq.sum.fubini.triangle2}) to
$a_{\left(  x,y\right)  }=1$. Thus, we obtain%
\[
\sum_{x=1}^{n}\ \ \sum_{y=x}^{n}1=\sum_{\left(  x,y\right)  \in Q_{n}}%
1=\sum_{y=1}^{n}\ \ \sum_{x=1}^{y}1.
\]
Hence,%
\[
\sum_{y=1}^{n}\ \ \sum_{x=1}^{y}1=\sum_{\left(  x,y\right)  \in Q_{n}%
}1=\left\vert Q_{n}\right\vert ,
\]
so that%
\[
\left\vert Q_{n}\right\vert =\sum_{y=1}^{n}\ \ \underbrace{\sum_{x=1}^{y}%
1}_{=y}=\sum_{y=1}^{n}y=\sum_{i=1}^{n}i=\dfrac{n\left(  n+1\right)  }%
{2}\ \ \ \ \ \ \ \ \ \ \left(  \text{by (\ref{eq.sum.littlegauss2})}\right)
.
\]

\end{itemize}

\textbf{Remarks:}

\begin{itemize}
\item The sum $\sum_{\left(  x,y\right)  \in Q_{n}}a_{\left(  x,y\right)  }$
in (\ref{eq.sum.fubini.triangle2}) can also be rewritten as $\sum
_{\substack{\left(  x,y\right)  \in\left\{  1,2,\ldots,n\right\}  ^{2};\\x\leq
y}}a_{\left(  x,y\right)  }$. It is also often written as $\sum_{1\leq x\leq
y\leq n}a_{\left(  x,y\right)  }$.

\item The proof of (\ref{eq.sum.fubini.triangle2}) is similar to that of
(\ref{eq.sum.fubini.triangle}).
\end{itemize}

\item \underline{\textbf{Fubini's theorem with a predicate:}} Let $X$ and $Y$
be two finite sets. For every pair $\left(  x,y\right)  \in X\times Y$, let
$\mathcal{A}\left(  x,y\right)  $ be a logical statement. For each $\left(
x,y\right)  \in X\times Y$ satisfying $\mathcal{A}\left(  x,y\right)  $, let
$a_{\left(  x,y\right)  }$ be an element of $\mathbb{A}$. Then,%
\begin{equation}
\sum_{x\in X}\ \ \sum_{\substack{y\in Y;\\\mathcal{A}\left(  x,y\right)
}}a_{\left(  x,y\right)  }=\sum_{\substack{\left(  x,y\right)  \in X\times
Y;\\\mathcal{A}\left(  x,y\right)  }}a_{\left(  x,y\right)  }=\sum_{y\in
Y}\ \ \sum_{\substack{x\in X;\\\mathcal{A}\left(  x,y\right)  }}a_{\left(
x,y\right)  }. \label{eq.sum.fubini.predicate}%
\end{equation}

\textbf{Examples:}

\begin{itemize}
\item For any $n\in\mathbb{N}$ and $m\in\mathbb{N}$, we have%
\begin{align*}
\sum_{x\in\left\{  1,2,\ldots,n\right\}  }\ \ \sum_{\substack{y\in\left\{
1,2,\ldots,m\right\}  ;\\x+y\text{ is even}}}xy  &  =\sum_{\substack{\left(
x,y\right)  \in\left\{  1,2,\ldots,n\right\}  \times\left\{  1,2,\ldots
,m\right\}  ;\\x+y\text{ is even}}}xy\\
&  =\sum_{y\in\left\{  1,2,\ldots,m\right\}  }\ \ \sum_{\substack{x\in\left\{
1,2,\ldots,n\right\}  ;\\x+y\text{ is even}}}xy.
\end{align*}
(This follows from (\ref{eq.sum.fubini.predicate}), applied to $X=\left\{
1,2,\ldots,n\right\}  $, $Y=\left\{  1,2,\ldots,m\right\}  $ and
$\mathcal{A}\left(  x,y\right)  =\left(  \text{\textquotedblleft}x+y\text{ is
even\textquotedblright}\right)  $.)
\end{itemize}

\textbf{Remarks:}

\begin{itemize}
\item We have assumed that the sets $X$ and $Y$ are finite. But
(\ref{eq.sum.fubini.predicate}) is still valid if we replace this assumption
by the weaker assumption that only finitely many $\left(  x,y\right)  \in
X\times Y$ satisfy $\mathcal{A}\left(  x,y\right)  $.

\item It is not hard to prove (\ref{eq.sum.fubini.predicate}) by suitably
adapting our proof of (\ref{eq.sum.fubini}).

\item The equality (\ref{eq.sum.fubini.triangle}) can be derived from
(\ref{eq.sum.fubini.predicate}) by setting $X=\left\{  1,2,\ldots,n\right\}
$, $Y=\left\{  1,2,\ldots,n\right\}  $ and $\mathcal{A}\left(  x,y\right)
=\left(  \text{\textquotedblleft}x+y\leq n\text{\textquotedblright}\right)  $.
Similarly, the equality (\ref{eq.sum.fubini.triangle2}) can be derived from
(\ref{eq.sum.fubini.predicate}) by setting $X=\left\{  1,2,\ldots,n\right\}
$, $Y=\left\{  1,2,\ldots,n\right\}  $ and $\mathcal{A}\left(  x,y\right)
=\left(  \text{\textquotedblleft}x\leq y\text{\textquotedblright}\right)  $.
\end{itemize}

\item \underline{\textbf{Interchange of predicates:}} Let $S$ be a finite set.
For every $s\in S$, let $\mathcal{A}\left(  s\right)  $ and $\mathcal{B}%
\left(  s\right)  $ be two equivalent logical statements. (\textquotedblleft
Equivalent\textquotedblright\ means that $\mathcal{A}\left(  s\right)  $ holds
if and only if $\mathcal{B}\left(  s\right)  $ holds.) Let $a_{s}$ be an
element of $\mathbb{A}$ for each $s\in S$. Then,%
\[
\sum_{\substack{s\in S;\\\mathcal{A}\left(  s\right)  }}a_{s}=\sum
_{\substack{s\in S;\\\mathcal{B}\left(  s\right)  }}a_{s}.
\]
(If you regard equivalent logical statements as identical, then you will see
this as a tautology. If not, it is still completely obvious, since the
equivalence of $\mathcal{A}\left(  s\right)  $ with $\mathcal{B}\left(
s\right)  $ shows that $\left\{  t\in S\ \mid\ \mathcal{A}\left(  t\right)
\right\}  =\left\{  t\in S\ \mid\ \mathcal{B}\left(  t\right)  \right\}  $.)

\item \underline{\textbf{Substituting the index I with a predicate:}} Let $S$
and $T$ be two finite sets. Let $f:S\rightarrow T$ be a \textbf{bijective}
map. Let $a_{t}$ be an element of $\mathbb{A}$ for each $t\in T$. For every
$t\in T$, let $\mathcal{A}\left(  t\right)  $ be a logical statement. Then,%
\begin{equation}
\sum_{\substack{t\in T;\\\mathcal{A}\left(  t\right)  }}a_{t}=\sum
_{\substack{s\in S;\\\mathcal{A}\left(  f\left(  s\right)  \right)
}}a_{f\left(  s\right)  }. \label{eq.sum.subs1-pred}%
\end{equation}

\textbf{Remarks:}

\begin{itemize}
\item The equality (\ref{eq.sum.subs1-pred}) is a generalization of
(\ref{eq.sum.subs1}). There is a similar generalization of (\ref{eq.sum.subs2}).

\item The equality (\ref{eq.sum.subs1-pred}) can be easily derived from
(\ref{eq.sum.subs1}). Indeed, let $S^{\prime}$ be the subset $\left\{  s\in
S\ \mid\ \mathcal{A}\left(  f\left(  s\right)  \right)  \right\}  $ of $S$,
and let $T^{\prime}$ be the subset $\left\{  t\in T\ \mid\ \mathcal{A}\left(
t\right)  \right\}  $ of $T$. Then, the map $S^{\prime}\rightarrow T^{\prime
},\ s\mapsto f\left(  s\right)  $ is well-defined and a
bijection\footnote{This is easy to see.}, and thus (\ref{eq.sum.subs1})
(applied to $S^{\prime}$, $T^{\prime}$ and this map instead of $S$, $T$ and
$f$) yields $\sum_{t\in T^{\prime}}a_{t}=\sum_{s\in S^{\prime}}a_{f\left(
s\right)  }$. But this is precisely the equality (\ref{eq.sum.subs1-pred}),
because clearly we have $\sum_{t\in T^{\prime}}=\sum_{\substack{t\in
T;\\\mathcal{A}\left(  t\right)  }}$ and $\sum_{s\in S^{\prime}}%
=\sum_{\substack{s\in S;\\\mathcal{A}\left(  f\left(  s\right)  \right)  }}$.
\end{itemize}
\end{itemize}

\subsubsection{Definition of $\prod$}

We shall now define the $\prod$ sign. Since the $\prod$ sign is (in many
aspects) analogous to the $\sum$ sign, we shall be brief and confine ourselves
to the bare necessities; we trust the reader to transfer most of what we said
about $\sum$ to the case of $\prod$. In particular, we shall give very few
examples and no proofs.

\begin{itemize}
\item If $S$ is a finite set, and if $a_{s}$ is an element of $\mathbb{A}$ for
each $s\in S$, then $\prod\nolimits_{s\in S}a_{s}$ denotes the product of all
of these elements $a_{s}$. Formally, this product is defined by recursion on
$\left\vert S\right\vert $, as follows:

\begin{itemize}
\item If $\left\vert S\right\vert =0$, then $\prod_{s\in S}a_{s}$ is defined
to be $1$.

\item Let $n\in\mathbb{N}$. Assume that we have defined $\prod_{s\in S}a_{s}$
for every finite set $S$ with $\left\vert S\right\vert =n$ (and every choice
of elements $a_{s}$ of $\mathbb{A}$). Now, if $S$ is a finite set with
$\left\vert S\right\vert =n+1$ (and if $a_{s}\in\mathbb{A}$ are chosen for all
$s\in S$), then $\prod_{s\in S}a_{s}$ is defined by picking any $t\in S$ and
setting%
\begin{equation}
\prod_{s\in S}a_{s}=a_{t}\cdot\prod_{s\in S\setminus\left\{  t\right\}  }%
a_{s}. \label{eq.prod.def.1}%
\end{equation}
As for $\sum_{s\in S}a_{s}$, this definition is not obviously legitimate, but
it can be proven to be legitimate nevertheless. (The proof is analogous to the
proof for $\sum_{s\in S}a_{s}$; see Subsection \ref{subsect.ind.gen-com.prods}
for details.)
\end{itemize}

\textbf{Examples:}

\begin{itemize}
\item If $S=\left\{  1,2,\ldots,n\right\}  $ (for some $n\in\mathbb{N}$) and
$a_{s}=s$ for every $s\in S$, then $\prod_{s\in S}a_{s}=\prod_{s\in S}%
s=1\cdot2\cdot\cdots\cdot n$. This number $1\cdot2\cdot\cdots\cdot n$ is
denoted by $n!$ and called the \textit{factorial of }$n$.

In particular,%
\begin{align*}
0!  &  =\prod_{s\in\left\{  1,2,\ldots,0\right\}  }s=\prod_{s\in\varnothing
}s\ \ \ \ \ \ \ \ \ \ \left(  \text{since }\left\{  1,2,\ldots,0\right\}
=\varnothing\right) \\
&  =1\ \ \ \ \ \ \ \ \ \ \left(  \text{since }\left\vert \varnothing
\right\vert =0\right)  ;\\
1!  &  =\prod_{s\in\left\{  1,2,\ldots,1\right\}  }s=\prod_{s\in\left\{
1\right\}  }s=1;\\
2!  &  =\prod_{s\in\left\{  1,2,\ldots,2\right\}  }s=\prod_{s\in\left\{
1,2\right\}  }s=1\cdot2=2;\\
3!  &  =\prod_{s\in\left\{  1,2,\ldots,3\right\}  }s=\prod_{s\in\left\{
1,2,3\right\}  }s=1\cdot2\cdot3=6;
\end{align*}
similarly,%
\[
4!=1\cdot2\cdot3\cdot
4=24;\ \ \ \ \ \ \ \ \ \ 5!=120;\ \ \ \ \ \ \ \ \ \ 6!=720;\ \ \ \ \ \ \ \ \ \ 7!=5040.
\]

Notice that%
\begin{equation}
n!=n\cdot\left(  n-1\right)  !\ \ \ \ \ \ \ \ \ \ \text{for any positive
integer }n. \label{eq.n!.rec}%
\end{equation}
(This can be obtained from (\ref{eq.prod.def.1}), applied to $S=\left\{
1,2,\ldots,n\right\}  $, $a_{s}=s$ and $t=n$.)
\end{itemize}

\textbf{Remarks:}

\begin{itemize}
\item The product $\prod_{s\in S}a_{s}$ is usually pronounced
\textquotedblleft product of the $a_{s}$ over all $s\in S$\textquotedblright%
\ or \textquotedblleft product of the $a_{s}$ with $s$ ranging over
$S$\textquotedblright\ or \textquotedblleft product of the $a_{s}$ with $s$
running through all elements of $S$\textquotedblright. The letter
\textquotedblleft$s$\textquotedblright\ in the product is called the
\textquotedblleft product index\textquotedblright, and its exact choice is
immaterial, as long as it does not already have a different meaning outside of
the product. The sign $\prod$ itself is called \textquotedblleft the product
sign\textquotedblright\ or \textquotedblleft the $\prod$
sign\textquotedblright. The numbers $a_{s}$ are called the \textit{factors} of
the product $\prod_{s\in S}a_{s}$. More precisely, for any given $t\in S$, we
can refer to the number $a_{t}$ as the \textquotedblleft factor corresponding
to the index $t$\textquotedblright\ (or as the \textquotedblleft factor for
$s=t$\textquotedblright, or as the \textquotedblleft factor for $t$%
\textquotedblright) of the product $\prod_{s\in S}a_{s}$.

\item When the set $S$ is empty, the product $\prod_{s\in S}a_{s}$ is called
an \textit{empty product}. Our definition implies that any empty product is
$1$. This convention is used throughout mathematics, except in rare occasions
where a slightly subtler version of it is used\footnote{Just as with sums, the
subtlety lies in the fact that mathematicians sometimes want an empty product
to be not the integer $1$ but the unity of some ring. As before, this does not
matter for us right now.}.

\item If $a\in\mathbb{A}$ and $n\in\mathbb{N}$, then the $n$-th power of $a$
(written $a^{n}$) is defined by%
\[
a^{n}=\underbrace{a\cdot a\cdot\cdots\cdot a}_{n\text{ times}}=\prod
_{i\in\left\{  1,2,\ldots,n\right\}  }a.
\]
Thus, $a^{0}$ is an empty product, and therefore equal to $1$. This holds for
any $a\in\mathbb{A}$, including $0$; thus, $0^{0}=1$. \textbf{There is nothing
controversial about the equality }$0^{0}=1$; it is a consequence of the only
reasonable definition of the $n$-th power of a number. Ignore anyone who tells
you that $0^{0}$ is \textquotedblleft undefined\textquotedblright\ or
\textquotedblleft indeterminate\textquotedblright\ or \textquotedblleft can be
$0$ or $1$ or anything, depending on the context\textquotedblright.\footnote{I
am talking about the \textbf{number} $0^{0}$ here. There is also something
called \textquotedblleft the
\href{https://en.wikipedia.org/wiki/Indeterminate form}{indeterminate form}
$0^{0}$\textquotedblright, which is a much different story.}

\item The product index (just like a summation index) needs not be a single
letter; it can be a pair or a triple, for example.

\item Mathematicians don't seem to have reached an agreement on the operator
precedence of the $\prod$ sign. My convention is that the product sign has
higher precedence than the plus sign (so an expression like $\prod_{s\in
S}a_{s}+b$ must be read as $\left(  \prod_{s\in S}a_{s}\right)  +b$, and not
as $\prod_{s\in S}\left(  a_{s}+b\right)  $); this is, of course, in line with
the standard convention that multiplication-like operations have higher
precedence than addition-like operations (\textquotedblleft
PEMDAS\textquotedblright). Be warned that some authors disagree even with this
convention. I strongly advise against writing things like $\prod_{s\in S}%
a_{s}b$, since it might mean both $\left(  \prod_{s\in S}a_{s}\right)  b$ and
$\prod_{s\in S}\left(  a_{s}b\right)  $ depending on the weather. In
particular, I advise against writing things like $\prod_{s\in S}a_{s}%
\cdot\prod_{s\in S}b_{s}$ without parentheses (although I do use a similar
convention for sums, namely $\sum_{s\in S}a_{s}+\sum_{s\in S}b_{s}$, and I
find it to be fairly harmless). These rules are not carved in stone, and you
should use whatever conventions make \textbf{you} safe from ambiguity; either
way, you should keep in mind that other authors make different choices.

\item An expression of the form \textquotedblleft$\prod_{s\in S}a_{s}%
$\textquotedblright\ (where $S$ is a finite set) is called a \textit{finite
product}.

\item We have required the set $S$ to be finite when defining $\prod_{s\in
S}a_{s}$. Such products are not generally defined when $S$ is infinite.
However, \textbf{some} infinite products can be made sense of. The simplest
case is when the set $S$ might be infinite, but only finitely many among the
$a_{s}$ are distinct from $1$. In this case, we can define $\prod_{s\in
S}a_{s}$ simply by discarding the factors which equal $1$ and multiplying the
finitely many remaining factors. Other situations in which infinite products
make sense appear in analysis and in topological algebra.

\item The product $\prod_{s\in S}a_{s}$ always belongs to $\mathbb{A}$.
\end{itemize}

\item A slightly more complicated version of the product sign is the
following: Let $S$ be a finite set, and let $\mathcal{A}\left(  s\right)  $ be
a logical statement defined for every $s\in S$. For each $s\in S$ satisfying
$\mathcal{A}\left(  s\right)  $, let $a_{s}$ be an element of $\mathbb{A}$.
Then, the product $\prod_{\substack{s\in S;\\\mathcal{A}\left(  s\right)
}}a_{s}$ is defined by%
\[
\prod_{\substack{s\in S;\\\mathcal{A}\left(  s\right)  }}a_{s}=\prod
_{s\in\left\{  t\in S\ \mid\ \mathcal{A}\left(  t\right)  \right\}  }a_{s}.
\]

\item Finally, here is the simplest version of the product sign: Let $u$ and
$v$ be two integers. As before, we understand the set $\left\{  u,u+1,\ldots
,v\right\}  $ to be empty when $u>v$. Let $a_{s}$ be an element of
$\mathbb{A}$ for each $s\in\left\{  u,u+1,\ldots,v\right\}  $. Then,
$\prod_{s=u}^{v}a_{s}$ is defined by%
\[
\prod_{s=u}^{v}a_{s}=\prod_{s\in\left\{  u,u+1,\ldots,v\right\}  }a_{s}.
\]

\textbf{Examples:}

\begin{itemize}
\item We have $\prod_{s=1}^{n}s=1\cdot2\cdot\cdots\cdot n=n!$ for each
$n\in\mathbb{N}$.
\end{itemize}

\textbf{Remarks:}

\begin{itemize}
\item The product $\prod_{s=u}^{v}a_{s}$ is usually pronounced
\textquotedblleft product of the $a_{s}$ for all $s$ from $u$ to $v$
(inclusive)\textquotedblright. It is often written $a_{u}\cdot a_{u+1}%
\cdot\cdots\cdot a_{v}$ (or just $a_{u}a_{u+1}\cdots a_{v}$), but this latter
notation has the same drawbacks as the similar notation $a_{u}+a_{u+1}%
+\cdots+a_{v}$ for $\sum_{s=u}^{v}a_{s}$.

\item The product $\prod_{s=u}^{v}a_{s}$ is said to be \textit{empty} whenever
$u>v$. As with sums, it does not matter how much smaller $v$ is than $u$; as
long as $v$ is smaller than $u$, the product is empty and equals $1$.
\end{itemize}
\end{itemize}

Thus we have introduced the main three forms of the product sign.

\subsubsection{Properties of $\prod$}

Now, let me summarize the most important properties of the $\prod$ sign. These
properties mirror the properties of $\sum$ discussed before; thus, I will
again be brief.

\begin{itemize}
\item \underline{\textbf{Splitting-off:}} Let $S$ be a finite set. Let $t\in
S$. Let $a_{s}$ be an element of $\mathbb{A}$ for each $s\in S$. Then,%
\[
\prod_{s\in S}a_{s}=a_{t}\cdot\prod_{s\in S\setminus\left\{  t\right\}  }%
a_{s}.
\]

\item \underline{\textbf{Splitting:}} Let $S$ be a finite set. Let $X$ and $Y
$ be two subsets of $S$ such that $X\cap Y=\varnothing$ and $X\cup Y=S$.
(Equivalently, $X$ and $Y$ are two subsets of $S$ such that each element of
$S$ lies in \textbf{exactly} one of $X$ and $Y$.) Let $a_{s}$ be an element of
$\mathbb{A}$ for each $s\in S$. Then,%
\[
\prod_{s\in S}a_{s}=\left(  \prod_{s\in X}a_{s}\right)  \cdot\left(
\prod_{s\in Y}a_{s}\right)  .
\]

\item \underline{\textbf{Splitting using a predicate:}} Let $S$ be a finite
set. Let $\mathcal{A}\left(  s\right)  $ be a logical statement for each $s\in
S$. Let $a_{s}$ be an element of $\mathbb{A}$ for each $s\in S$. Then,%
\[
\prod_{s\in S}a_{s}=\left(  \prod_{\substack{s\in S; \\\mathcal{A}\left(
s\right)  }}a_{s}\right)  \cdot\left(  \prod_{\substack{s\in S; \\\text{not
}\mathcal{A}\left(  s\right)  }}a_{s}\right)  .
\]

\item \underline{\textbf{Multiplying equal values:}} Let $S$ be a finite set.
Let $a$ be an element of $\mathbb{A}$. Then,%
\[
\prod_{s\in S}a=a^{\left\vert S\right\vert }.
\]

\item \underline{\textbf{Splitting a factor:}} Let $S$ be a finite set. For
every $s\in S$, let $a_{s}$ and $b_{s}$ be elements of $\mathbb{A}$. Then,%
\begin{equation}
\prod_{s\in S}\left(  a_{s}b_{s}\right)  =\left(  \prod_{s\in S}a_{s}\right)
\cdot\left(  \prod_{s\in S}b_{s}\right)  . \label{eq.prod.linear1}%
\end{equation}

\textbf{Examples:}

\begin{itemize}
\item Here is a frequently used particular case of (\ref{eq.prod.linear1}):
Let $S$ be a finite set. For every $s\in S$, let $b_{s}$ be an element of
$\mathbb{A}$. Let $a$ be an element of $\mathbb{A}$. Then,
(\ref{eq.prod.linear1}) (applied to $a_{s}=a$) yields%
\begin{equation}
\prod_{s\in S}\left(  ab_{s}\right)  =\underbrace{\left(  \prod_{s\in
S}a\right)  }_{=a^{\left\vert S\right\vert }}\cdot\left(  \prod_{s\in S}%
b_{s}\right)  =a^{\left\vert S\right\vert }\cdot\left(  \prod_{s\in S}%
b_{s}\right)  . \label{eq.prod.linear1.ex1}%
\end{equation}

\item Here is an even further particular case: Let $S$ be a finite set. For
every $s\in S$, let $b_{s}$ be an element of $\mathbb{A}$. Then,%
\[
\prod_{s\in S}\underbrace{\left(  -b_{s}\right)  }_{=\left(  -1\right)  b_{s}%
}=\prod_{s\in S}\left(  \left(  -1\right)  b_{s}\right)  =\left(  -1\right)
^{\left\vert S\right\vert }\cdot\left(  \prod_{s\in S}b_{s}\right)
\]
(by (\ref{eq.prod.linear1.ex1}), applied to $a=-1$).
\end{itemize}

\item \underline{\textbf{Factoring out an exponent:}} Let $S$ be a finite set.
For every $s\in S$, let $a_{s}$ be an element of $\mathbb{A}$. Also, let
$\lambda\in\mathbb{N}$. Then,%
\[
\prod_{s\in S}a_{s}^{\lambda}=\left(  \prod_{s\in S}a_{s}\right)  ^{\lambda}.
\]

\item \underline{\textbf{Factoring out an integer exponent:}} Let $S$ be a
finite set. For every $s\in S$, let $a_{s}$ be a nonzero element of
$\mathbb{A}$. Also, let $\lambda\in\mathbb{Z}$. Then,%
\[
\prod_{s\in S}a_{s}^{\lambda}=\left(  \prod_{s\in S}a_{s}\right)  ^{\lambda}.
\]

\textbf{Remark:} Applying this to $\lambda=-1$, we obtain%
\[
\prod_{s\in S}a_{s}^{-1}=\left(  \prod_{s\in S}a_{s}\right)  ^{-1}.
\]
In other words,%
\[
\prod_{s\in S}\dfrac{1}{a_{s}}=\dfrac{1}{\prod_{s\in S}a_{s}}.
\]

\item \underline{\textbf{Ones multiply to one:}} Let $S$ be a finite set.
Then,%
\[
\prod_{s\in S}1=1.
\]

\item \underline{\textbf{Dropping ones:}} Let $S$ be a finite set. Let $a_{s}$
be an element of $\mathbb{A}$ for each $s\in S$. Let $T$ be a subset of $S$
such that every $s\in T$ satisfies $a_{s}=1$. Then,%
\[
\prod_{s\in S}a_{s}=\prod_{s\in S\setminus T}a_{s}.
\]

\item \underline{\textbf{Renaming the index:}} Let $S$ be a finite set. Let
$a_{s}$ be an element of $\mathbb{A}$ for each $s\in S$. Then,%
\[
\prod_{s\in S}a_{s}=\prod_{t\in S}a_{t}.
\]

\item \underline{\textbf{Substituting the index I:}} Let $S$ and $T$ be two
finite sets. Let $f:S\rightarrow T$ be a \textbf{bijective} map. Let $a_{t}$
be an element of $\mathbb{A}$ for each $t\in T$. Then,%
\[
\prod_{t\in T}a_{t}=\prod_{s\in S}a_{f\left(  s\right)  }.
\]

\item \underline{\textbf{Substituting the index II:}} Let $S$ and $T$ be two
finite sets. Let $f:S\rightarrow T$ be a \textbf{bijective} map. Let $a_{s}$
be an element of $\mathbb{A}$ for each $s\in S$. Then,%
\[
\prod_{s\in S}a_{s}=\prod_{t\in T}a_{f^{-1}\left(  t\right)  }.
\]

\item \underline{\textbf{Telescoping products:}} Let $u$ and $v$ be two
integers such that $u-1\leq v$. Let $a_{s}$ be an element of $\mathbb{A}$ for
each $s\in\left\{  u-1,u,\ldots,v\right\}  $. Then,%
\begin{equation}
\prod_{s=u}^{v}\dfrac{a_{s}}{a_{s-1}}=\dfrac{a_{v}}{a_{u-1}}
\label{eq.prod.telescope}%
\end{equation}
(provided that $a_{s-1}\neq0$ for all $s\in\left\{  u,u+1,\ldots,v\right\}  $).

\textbf{Examples:}

\begin{itemize}
\item Let $n$ be a positive integer. Then,%
\begin{align*}
\prod_{s=2}^{n}\underbrace{\left(  1-\dfrac{1}{s}\right)  }_{=\dfrac{s-1}%
{s}=\dfrac{1/s}{1/\left(  s-1\right)  }}  &  =\prod_{s=2}^{n}\dfrac
{1/s}{1/\left(  s-1\right)  }=\dfrac{1/n}{1/\left(  2-1\right)  }\\
&  \ \ \ \ \ \ \ \ \ \ \left(  \text{by (\ref{eq.prod.telescope}), applied to
}u=2\text{, }v=n\text{ and }a_{s}=1/s\right) \\
&  =\dfrac{1}{n}.
\end{align*}

\end{itemize}

\item \underline{\textbf{Restricting to a subset:}} Let $S$ be a finite set.
Let $T$ be a subset of $S$. Let $a_{s}$ be an element of $\mathbb{A}$ for each
$s\in T$. Then,%
\[
\prod_{\substack{s\in S;\\s\in T}}a_{s}=\prod_{s\in T}a_{s}.
\]
\textbf{Remark:} Here is a slightly more general form of this rule: Let $S$ be
a finite set. Let $T$ be a subset of $S$. Let $\mathcal{A}\left(  s\right)  $
be a logical statement for each $s\in S$. Let $a_{s}$ be an element of
$\mathbb{A}$ for each $s\in T$ satisfying $\mathcal{A}\left(  s\right)  $.
Then,%
\[
\prod_{\substack{s\in S;\\s\in T;\\\mathcal{A}\left(  s\right)  }}a_{s}%
=\prod_{\substack{s\in T;\\\mathcal{A}\left(  s\right)  }}a_{s}.
\]

\item \underline{\textbf{Splitting a product by a value of a function:}} Let
$S$ be a finite set. Let $W$ be a set. Let $f:S\rightarrow W$ be a map. Let
$a_{s}$ be an element of $\mathbb{A}$ for each $s\in S$. Then,%
\[
\prod_{s\in S}a_{s}=\prod_{w\in W}\ \ \prod_{\substack{s\in S;\\f\left(
s\right)  =w}}a_{s}.
\]
(The right hand side is to be read as $\prod_{w\in W}\left(  \prod
_{\substack{s\in S;\\f\left(  s\right)  =w}}a_{s}\right)  $.)

\item \underline{\textbf{Splitting a product into subproducts:}} Let $S$ be a
finite set. Let $S_{1},S_{2},\ldots,S_{n}$ be finitely many subsets of $S$.
Assume that these subsets $S_{1},S_{2},\ldots,S_{n}$ are pairwise disjoint
(i.e., we have $S_{i}\cap S_{j}=\varnothing$ for any two distinct elements $i$
and $j$ of $\left\{  1,2,\ldots,n\right\}  $) and their union is $S$. (Thus,
every element of $S$ lies in precisely one of the subsets $S_{1},S_{2}%
,\ldots,S_{n}$.) Let $a_{s}$ be an element of $\mathbb{A}$ for each $s\in S$.
Then,%
\[
\prod_{s\in S}a_{s}=\prod_{w=1}^{n}\ \ \prod_{s\in S_{w}}a_{s}.
\]

\item \underline{\textbf{Fubini's theorem (interchanging the order of
multiplication):}} Let $X$ and $Y$ be two finite sets. Let $a_{\left(
x,y\right)  }$ be an element of $\mathbb{A}$ for each $\left(  x,y\right)  \in
X\times Y$. Then,%
\[
\prod_{x\in X}\ \ \prod_{y\in Y}a_{\left(  x,y\right)  }=\prod_{\left(
x,y\right)  \in X\times Y}a_{\left(  x,y\right)  }=\prod_{y\in Y}%
\ \ \prod_{x\in X}a_{\left(  x,y\right)  }.
\]

In particular, if $n$ and $m$ are two elements of $\mathbb{N}$, and if
$a_{\left(  x,y\right)  }$ is an element of $\mathbb{A}$ for each $\left(
x,y\right)  \in\left\{  1,2,\ldots,n\right\}  \times\left\{  1,2,\ldots
,m\right\}  $, then%
\[
\prod_{x=1}^{n}\ \ \prod_{y=1}^{m}a_{\left(  x,y\right)  }=\prod_{\left(
x,y\right)  \in\left\{  1,2,\ldots,n\right\}  \times\left\{  1,2,\ldots
,m\right\}  }a_{\left(  x,y\right)  }=\prod_{y=1}^{m}\ \ \prod_{x=1}%
^{n}a_{\left(  x,y\right)  }.
\]

\item \underline{\textbf{Triangular Fubini's theorem I:}} Let $n\in\mathbb{N}%
$. Let $T_{n}$ be the set \newline$\left\{  \left(  x,y\right)  \in\left\{
1,2,3,\ldots\right\}  ^{2}\ \mid\ x+y\leq n\right\}  $. Let $a_{\left(
x,y\right)  }$ be an element of $\mathbb{A}$ for each $\left(  x,y\right)  \in
T_{n}$. Then,%
\[
\prod_{x=1}^{n}\ \ \prod_{y=1}^{n-x}a_{\left(  x,y\right)  }=\prod_{\left(
x,y\right)  \in T_{n}}a_{\left(  x,y\right)  }=\prod_{y=1}^{n}\ \ \prod
_{x=1}^{n-y}a_{\left(  x,y\right)  }.
\]

\item \underline{\textbf{Triangular Fubini's theorem II:}} Let $n\in
\mathbb{N}$. Let $Q_{n}$ be the set \newline$\left\{  \left(  x,y\right)
\in\left\{  1,2,\ldots,n\right\}  ^{2}\ \mid\ x\leq y\right\}  $. Let
$a_{\left(  x,y\right)  }$ be an element of $\mathbb{A}$ for each $\left(
x,y\right)  \in Q_{n}$. Then,%
\[
\prod_{x=1}^{n}\ \ \prod_{y=x}^{n}a_{\left(  x,y\right)  }=\prod_{\left(
x,y\right)  \in Q_{n}}a_{\left(  x,y\right)  }=\prod_{y=1}^{n}\ \ \prod
_{x=1}^{y}a_{\left(  x,y\right)  }.
\]

\item \underline{\textbf{Fubini's theorem with a predicate:}} Let $X$ and $Y$
be two finite sets. For every pair $\left(  x,y\right)  \in X\times Y$, let
$\mathcal{A}\left(  x,y\right)  $ be a logical statement. For each $\left(
x,y\right)  \in X\times Y$ satisfying $\mathcal{A}\left(  x,y\right)  $, let
$a_{\left(  x,y\right)  }$ be an element of $\mathbb{A}$. Then,%
\[
\prod_{x\in X}\ \ \prod_{\substack{y\in Y;\\\mathcal{A}\left(  x,y\right)
}}a_{\left(  x,y\right)  }=\prod_{\substack{\left(  x,y\right)  \in X\times
Y;\\\mathcal{A}\left(  x,y\right)  }}a_{\left(  x,y\right)  }=\prod_{y\in
Y}\ \ \prod_{\substack{x\in X;\\\mathcal{A}\left(  x,y\right)  }}a_{\left(
x,y\right)  }.
\]

\item \underline{\textbf{Interchange of predicates:}} Let $S$ be a finite set.
For every $s\in S$, let $\mathcal{A}\left(  s\right)  $ and $\mathcal{B}%
\left(  s\right)  $ be two equivalent logical statements. (\textquotedblleft
Equivalent\textquotedblright\ means that $\mathcal{A}\left(  s\right)  $ holds
if and only if $\mathcal{B}\left(  s\right)  $ holds.) Let $a_{s}$ be an
element of $\mathbb{A}$ for each $s\in S$. Then,%
\[
\prod_{\substack{s\in S;\\\mathcal{A}\left(  s\right)  }}a_{s}=\prod
_{\substack{s\in S;\\\mathcal{B}\left(  s\right)  }}a_{s}.
\]

\item \underline{\textbf{Substituting the index I with a predicate:}} Let $S$
and $T$ be two finite sets. Let $f:S\rightarrow T$ be a \textbf{bijective}
map. Let $a_{t}$ be an element of $\mathbb{A}$ for each $t\in T$. For every
$t\in T$, let $\mathcal{A}\left(  t\right)  $ be a logical statement. Then,%
\[
\prod_{\substack{t\in T;\\\mathcal{A}\left(  t\right)  }}a_{t}=\prod
_{\substack{s\in S;\\\mathcal{A}\left(  f\left(  s\right)  \right)
}}a_{f\left(  s\right)  }.
\]

\end{itemize}

\subsection{\label{sect.polynomials-emergency}Polynomials: a precise
definition}

As I have already mentioned in the above list of prerequisites, the notion of
polynomials (in one and in several indeterminates) will be occasionally used
in these notes. Most likely, the reader already has at least a vague
understanding of this notion (e.g., from high school); this vague
understanding is probably sufficient for reading most of these notes. But
polynomials are one of the most important notions in algebra (if not to say in
mathematics), and the reader will likely encounter them over and over; sooner
or later, it will happen that the vague understanding is not sufficient and
some subtleties do matter. For that reason, anyone serious about doing
abstract algebra should know a complete and correct definition of polynomials
and have some experience working with it. I shall not give a complete
definition of the most general notion of polynomials in these notes, but I
will comment on some of the subtleties and define an important special case
(that of polynomials in one variable with rational coefficients) in the
present section. A reader is probably best advised to skip this section on
their first read.

It is not easy to find a good (formal and sufficiently general) treatment of
polynomials in textbooks. Various authors tend to skimp on subtleties and
technical points such as the notion of an \textquotedblleft
indeterminate\textquotedblright, or the precise meaning of \textquotedblleft
formal expression\textquotedblright\ in the slogan \textquotedblleft a
polynomial is a formal expression\textquotedblright\ (the best texts do not
use this vague slogan at all), or the definition of the degree of the zero
polynomial, or the difference between regarding polynomials as sequences
(which is the classical viewpoint and particularly useful for polynomials in
one variable) and regarding polynomials as elements of a monoid ring (which is
important in the case of several variables, since it allows us to regard the
polynomial rings $\mathbb{Q}\left[  X\right]  $ and $\mathbb{Q}\left[
Y\right]  $ as two distinct subrings of $\mathbb{Q}\left[  X,Y\right]  $).
They also tend to take some questionable shortcuts, such as defining
polynomials in $n$ variables (by induction over $n$) as one-variable
polynomials over the ring of $\left(  n-1\right)  $-variable polynomials (this
shortcut has several shortcomings, such as making the symmetric role of the
$n$ variables opaque, and functioning only for finitely many variables).

More often than not, the polynomials we will be using will be polynomials in
one variable. These are usually handled well in good books on abstract algebra
-- e.g., in \cite[\S 4.5]{Walker87}, in \cite[Appendix G]{Hungerford}, in
\cite[Chapter III, \S 5]{Hungerford-03}, in \cite[Chapter A-3]{Rotman15}, in
\cite[\S 4.1, \S 4.2]{HoffmanKunze} (although in \cite[\S 4.1, \S 4.2]%
{HoffmanKunze}, only polynomials over fields are studied, but the definition
applies to commutative rings mutatis mutandis), in \cite[\S 8]{AmaEsc05}, and
in \cite[Chapter III, \S 6]{BirkMac}. Most of these treatments rely on the
notion of a \textit{commutative ring}, which is not difficult but somewhat
abstract (I shall introduce it below in Section \ref{sect.commring}).

Let me give a brief survey of the notion of univariate polynomials (i.e.,
polynomials in one variable). I shall define them as sequences. For the sake
of simplicity, I shall only talk of polynomials with rational coefficients.
Similarly, one can define polynomials with integer coefficients, with real
coefficients, or with complex coefficients; of course, one then has to replace
each \textquotedblleft$\mathbb{Q}$\textquotedblright\ by a \textquotedblleft%
$\mathbb{Z}$\textquotedblright, an \textquotedblleft$\mathbb{R}$%
\textquotedblright\ or a \textquotedblleft$\mathbb{C}$\textquotedblright.

The rough idea behind the definition of a polynomial is that a polynomial with
rational coefficients should be a \textquotedblleft formal
expression\textquotedblright\ which is built out of rational numbers, an
\textquotedblleft indeterminate\textquotedblright\ $X$ as well as addition,
subtraction and multiplication signs, such as $X^{4}-27X+\dfrac{3}{2}$ or
$-X^{3}+2X+1$ or $\dfrac{1}{3}\left(  X-3\right)  \cdot X^{2}$ or
$X^{4}+7X^{3}\left(  X-2\right)  $ or $-15$. We have not explicitly allowed
powers, but we understand $X^{n}$ to mean the product $\underbrace{XX\cdots
X}_{n\text{ times}}$ (which is $1$ when $n=0$). Notice that division is not
allowed, so we cannot get $\dfrac{X}{X+1}$ (but we can get $\dfrac{3}{2}X$,
because $\dfrac{3}{2}$ is a rational number). Notice also that a polynomial
can be a single rational number, since we never said that $X$ must necessarily
be used; for instance, $-15$ and $0$ are polynomials.

This is, of course, not a valid definition. One problem with it that it does
not explain what a \textquotedblleft formal expression\textquotedblright\ is.
For starters, we want an expression that is well-defined -- i.e., into that we
can substitute a rational number for $X$ and obtain a valid term. For example,
$X-+\cdot5$ is not well-defined, so it does not fit our bill; neither is the
\textquotedblleft empty expression\textquotedblright. Furthermore, when do we
want two \textquotedblleft formal expressions\textquotedblright\ to be viewed
as one and the same polynomial? Do we want to equate $X\left(  X+2\right)  $
with $X^{2}+2X$ ? Do we want to equate $0X^{3}+2X+1$ with $2X+1$ ? The answer
is \textquotedblleft yes\textquotedblright\ both times, but a general rule is
not easy to give if we keep talking of \textquotedblleft formal
expressions\textquotedblright.

We \textit{could} define two polynomials $p\left(  X\right)  $ and $q\left(
X\right)  $ to be equal if and only if, for every number $\alpha\in\mathbb{Q}%
$, the values $p\left(  \alpha\right)  $ and $q\left(  \alpha\right)  $
(obtained by substituting $\alpha$ for $X$ in $p$ and in $q$, respectively)
are equal. This would be tantamount to treating polynomials as
\textit{functions}: it would mean that we identify a polynomial $p\left(
X\right)  $ with the function $\mathbb{Q}\rightarrow\mathbb{Q},\ \alpha\mapsto
p\left(  \alpha\right)  $. Such a definition would work well as long as we
would do only rather basic things with it\footnote{And some authors, such as
Axler in \cite[Chapter 4]{Axler}, do use this definition.}, but as soon as we
would try to go deeper, we would encounter technical issues which would make
it inadequate and painful\footnote{Here are three of these issues:
\par
\begin{itemize}
\item One of the strengths of polynomials is that we can evaluate them not
only at numbers, but also at many other things, e.g., at square matrices:
Evaluating the polynomial $X^{2}-3X$ at the square matrix $\left(
\begin{array}
[c]{cc}%
1 & 3\\
-1 & 2
\end{array}
\right)  $ gives $\left(
\begin{array}
[c]{cc}%
1 & 3\\
-1 & 2
\end{array}
\right)  ^{2}-3\left(
\begin{array}
[c]{cc}%
1 & 3\\
-1 & 2
\end{array}
\right)  =\left(
\begin{array}
[c]{cc}%
-5 & 0\\
0 & -5
\end{array}
\right)  $. However, a function must have a well-defined domain, and does not
make sense outside of this domain. So, if the polynomial $X^{2}-3X$ is
regarded as the function $\mathbb{Q}\rightarrow\mathbb{Q},\ \alpha
\mapsto\alpha^{2}-3\alpha$, then it makes no sense to evaluate this polynomial
at the matrix $\left(
\begin{array}
[c]{cc}%
1 & 3\\
-1 & 2
\end{array}
\right)  $, just because this matrix does not lie in the domain $\mathbb{Q}$
of the function. We could, of course, extend the domain of the function to
(say) the set of square matrices over $\mathbb{Q}$, but then we would still
have the same problem with other things that we want to evaluate polynomials
at. At some point we want to be able to evaluate polynomials at functions and
at other polynomials, and if we would try to achieve this by extending the
domain, we would have to do this over and over, because each time we extend
the domain, we get even more polynomials to evaluate our polynomials at; thus,
the definition would be eternally \textquotedblleft hunting its own
tail\textquotedblright! (We could resolve this difficulty by defining
polynomials as \textit{natural transformations} in the sense of category
theory. I do not want to even go into this definition here, as it would take
several pages to properly introduce. At this point, it is not worth the
hassle.)
\par
\item Let $p\left(  X\right)  $ be a polynomial with real coefficients. Then,
it should be obvious that $p\left(  X\right)  $ can also be viewed as a
polynomial with complex coefficients: For instance, if $p\left(  X\right)  $
was defined as $3X+\dfrac{7}{2}X\left(  X-1\right)  $, then we can view the
numbers $3$, $\dfrac{7}{2}$ and $-1$ appearing in its definition as complex
numbers, and thus get a polynomial with complex coefficients. But wait! What
if two polynomials $p\left(  X\right)  $ and $q\left(  X\right)  $ are equal
when viewed as polynomials with real coefficients, but become distinct when
viewed as polynomials with complex coefficients (because when we view them as
polynomials with complex coefficients, their domains grow larger to include
complex numbers, and a new complex $\alpha$ might perhaps no longer satisfy
$p\left(  \alpha\right)  =q\left(  \alpha\right)  $ )? This does not actually
happen, but ruling this out is not obvious if you regard polynomials as
functions.
\par
\item (This requires some familiarity with finite fields:) Treating
polynomials as functions works reasonably well for polynomials with integer,
rational, real and complex coefficients (as long as one is not too demanding).
But we will eventually want to consider polynomials with coefficients in any
arbitrary commutative ring $\mathbb{K}$. An example for a commutative ring
$\mathbb{K}$ is the finite field $\mathbb{F}_{p}$ with $p$ elements, where $p$
is a prime. (This finite field $\mathbb{F}_{p}$ is better known as the ring of
integers modulo $p$.) If we define polynomials with coefficients in
$\mathbb{F}_{p}$ as functions $\mathbb{F}_{p}\rightarrow\mathbb{F}_{p}$, then
we really run into problems; for example, the polynomials $X$ and $X^{p}$ over
this field become identical as functions!
\end{itemize}
}. Also, if we equated polynomials with the functions they describe, then we
would waste the word \textquotedblleft polynomial\textquotedblright\ on a
concept (a function described by a polynomial) that already has a word for it
(namely, \textit{polynomial function}).

The preceding paragraphs indicate that it is worth defining \textquotedblleft
polynomials\textquotedblright\ in a way that, on the one hand, conveys the
idea that they are more \textquotedblleft formal expressions\textquotedblright%
\ than \textquotedblleft functions\textquotedblright, but on the other hand,
is less nebulous than \textquotedblleft formal expression\textquotedblright.
Here is one such definition:

\begin{definition}
\label{def.polynomial-univar}\textbf{(a)} A \textit{univariate polynomial with
rational coefficients} means a sequence $\left(  p_{0},p_{1},p_{2}%
,\ldots\right)  \in\mathbb{Q}^{\infty}$ of elements of $\mathbb{Q}$ such that%
\begin{equation}
\text{all but finitely many }k\in\mathbb{N}\text{ satisfy }p_{k}=0.
\label{eq.def.polynomial-univar.finite}%
\end{equation}
Here, the phrase \textquotedblleft all but finitely many $k\in\mathbb{N}$
satisfy $p_{k}=0$\textquotedblright\ means \textquotedblleft there exists some
finite subset $J$ of $\mathbb{N}$ such that every $k\in\mathbb{N}\setminus J$
satisfies $p_{k}=0$\textquotedblright. (See Definition \ref{def.allbutfin} for
the general definition of \textquotedblleft all but finitely
many\textquotedblright, and Section \ref{sect.infperm} for some practice with
this concept.) More concretely, the condition
(\ref{eq.def.polynomial-univar.finite}) can be rewritten as follows: The
sequence $\left(  p_{0},p_{1},p_{2},\ldots\right)  $ contains only zeroes from
some point on (i.e., there exists some $N\in\mathbb{N}$ such that
$p_{N}=p_{N+1}=p_{N+2}=\cdots=0$).

For the remainder of this definition, \textquotedblleft univariate polynomial
with rational coefficients\textquotedblright\ will be abbreviated as
\textquotedblleft polynomial\textquotedblright.

For example, the sequences $\left(  0,0,0,\ldots\right)  $, $\left(
1,3,5,0,0,0,\ldots\right)  $, $\left(  4,0,-\dfrac{2}{3},5,0,0,0,\ldots
\right)  $, $\left(  0,-1,\dfrac{1}{2},0,0,0,\ldots\right)  $ (where the
\textquotedblleft$\ldots$\textquotedblright\ stand for infinitely many zeroes)
are polynomials, but the sequence $\left(  1,1,1,\ldots\right)  $ (where the
\textquotedblleft$\ldots$\textquotedblright\ stands for infinitely many $1$'s)
is not (since it does not satisfy (\ref{eq.def.polynomial-univar.finite})).

So we have defined a polynomial as an infinite sequence of rational numbers
with a certain property. So far, this does not seem to reflect any intuition
of polynomials as \textquotedblleft formal expressions\textquotedblright.
However, we shall soon (namely, in Definition \ref{def.polynomial-univar}
\textbf{(j)}) identify the polynomial $\left(  p_{0},p_{1},p_{2}%
,\ldots\right)  \in\mathbb{Q}^{\infty}$ with the \textquotedblleft formal
expression\textquotedblright\ $p_{0}+p_{1}X+p_{2}X^{2}+\cdots$ (this is an
infinite sum, but due to (\ref{eq.def.polynomial-univar.finite}) all but its
first few terms are $0$ and thus can be neglected). For instance, the
polynomial $\left(  1,3,5,0,0,0,\ldots\right)  $ will be identified with the
\textquotedblleft formal expression\textquotedblright\ $1+3X+5X^{2}%
+0X^{3}+0X^{4}+0X^{5}+\cdots=1+3X+5X^{2}$. Of course, we cannot do this
identification right now, since we do not have a reasonable definition of $X$.

\textbf{(b)} We let $\mathbb{Q}\left[  X\right]  $ denote the set of all
univariate polynomials with rational coefficients. Given a polynomial
$p=\left(  p_{0},p_{1},p_{2},\ldots\right)  \in\mathbb{Q}\left[  X\right]  $,
we denote the numbers $p_{0},p_{1},p_{2},\ldots$ as the \textit{coefficients}
of $p$. More precisely, for every $i\in\mathbb{N}$, we shall refer to $p_{i}$
as the $i$\textit{-th coefficient} of $p$. (Do not forget that we are counting
from $0$ here: any polynomial \textquotedblleft begins\textquotedblright\ with
its $0$-th coefficient.) The $0$-th coefficient of $p$ is also known as the
\textit{constant term} of $p$.

Instead of \textquotedblleft the $i$-th coefficient of $p$\textquotedblright,
we often also say \textquotedblleft the \textit{coefficient before }$X^{i}%
$\textit{ of }$p$\textquotedblright\ or \textquotedblleft the
\textit{coefficient of }$X^{i}$ \textit{in }$p$\textquotedblright.

Thus, any polynomial $p\in\mathbb{Q}\left[  X\right]  $ is the sequence of its coefficients.

\textbf{(c)} We denote the polynomial $\left(  0,0,0,\ldots\right)
\in\mathbb{Q}\left[  X\right]  $ by $\mathbf{0}$. We will also write $0$ for
it when no confusion with the number $0$ is possible. The polynomial
$\mathbf{0}$ is called the \textit{zero polynomial}. A polynomial
$p\in\mathbb{Q}\left[  X\right]  $ is said to be \textit{nonzero} if
$p\neq\mathbf{0}$.

\textbf{(d)} We denote the polynomial $\left(  1,0,0,0,\ldots\right)
\in\mathbb{Q}\left[  X\right]  $ by $\mathbf{1}$. We will also write $1$ for
it when no confusion with the number $1$ is possible.

\textbf{(e)} For any $\lambda\in\mathbb{Q}$, we denote the polynomial $\left(
\lambda,0,0,0,\ldots\right)  \in\mathbb{Q}\left[  X\right]  $ by
$\operatorname*{const}\lambda$. We call it the \textit{constant polynomial
with value }$\lambda$. It is often useful to identify $\lambda\in\mathbb{Q}$
with $\operatorname*{const}\lambda\in\mathbb{Q}\left[  X\right]  $. Notice
that $\mathbf{0}=\operatorname*{const}0$ and $\mathbf{1}=\operatorname*{const}%
1$.

\textbf{(f)} Now, let us define the sum, the difference and the product of two
polynomials. Indeed, let $a=\left(  a_{0},a_{1},a_{2},\ldots\right)
\in\mathbb{Q}\left[  X\right]  $ and $b=\left(  b_{0},b_{1},b_{2}%
,\ldots\right)  \in\mathbb{Q}\left[  X\right]  $ be two polynomials. Then, we
define three polynomials $a+b$, $a-b$ and $a\cdot b$ in $\mathbb{Q}\left[
X\right]  $ by%
\begin{align*}
a+b  &  =\left(  a_{0}+b_{0},a_{1}+b_{1},a_{2}+b_{2},\ldots\right)  ;\\
a-b  &  =\left(  a_{0}-b_{0},a_{1}-b_{1},a_{2}-b_{2},\ldots\right)  ;\\
a\cdot b  &  =\left(  c_{0},c_{1},c_{2},\ldots\right)  ,
\end{align*}
where%
\[
c_{k}=\sum_{i=0}^{k}a_{i}b_{k-i}\ \ \ \ \ \ \ \ \ \ \text{for every }%
k\in\mathbb{N}.
\]
We call $a+b$ the \textit{sum} of $a$ and $b$; we call $a-b$ the
\textit{difference} of $a$ and $b$; we call $a\cdot b$ the \textit{product} of
$a$ and $b$. We abbreviate $a\cdot b$ by $ab$, and we abbreviate
$\mathbf{0}-a$ by $-a$.

For example,%
\begin{align*}
\left(  1,2,2,0,0,\ldots\right)  +\left(  3,0,-1,0,0,0,\ldots\right)   &
=\left(  4,2,1,0,0,0,\ldots\right)  ;\\
\left(  1,2,2,0,0,\ldots\right)  -\left(  3,0,-1,0,0,0,\ldots\right)   &
=\left(  -2,2,3,0,0,0,\ldots\right)  ;\\
\left(  1,2,2,0,0,\ldots\right)  \cdot\left(  3,0,-1,0,0,0,\ldots\right)   &
=\left(  3,6,5,-2,-2,0,0,0,\ldots\right)  .
\end{align*}

The definition of $a+b$ essentially says that \textquotedblleft polynomials
are added coefficientwise\textquotedblright\ (i.e., in order to obtain the sum
of two polynomials $a$ and $b$, it suffices to add each coefficient of $a$ to
the corresponding coefficient of $b$). Similarly, the definition of $a-b$ says
the same thing about subtraction. The definition of $a\cdot b$ is more
surprising. However, it loses its mystique when we identify the polynomials
$a$ and $b$ with the \textquotedblleft formal expressions\textquotedblright%
\ $a_{0}+a_{1}X+a_{2}X^{2}+\cdots$ and $b_{0}+b_{1}X+b_{2}X^{2}+\cdots$
(although, at this point, we do not know what these expressions really mean);
indeed, it simply says that
\[
\left(  a_{0}+a_{1}X+a_{2}X^{2}+\cdots\right)  \left(  b_{0}+b_{1}X+b_{2}%
X^{2}+\cdots\right)  =c_{0}+c_{1}X+c_{2}X^{2}+\cdots,
\]
where $c_{k}=\sum_{i=0}^{k}a_{i}b_{k-i}$ for every $k\in\mathbb{N}$. This is
precisely what one would expect, because if you expand $\left(  a_{0}%
+a_{1}X+a_{2}X^{2}+\cdots\right)  \left(  b_{0}+b_{1}X+b_{2}X^{2}%
+\cdots\right)  $ using the distributive law and collect equal powers of $X$,
then you get precisely $c_{0}+c_{1}X+c_{2}X^{2}+\cdots$. Thus, the definition
of $a\cdot b$ has been tailored to make the distributive law hold.

(By the way, why is $a\cdot b$ a polynomial? That is, why does it satisfy
(\ref{eq.def.polynomial-univar.finite}) ? The proof is easy, but we omit it.)

Addition, subtraction and multiplication of polynomials satisfy some of the
same rules as addition, subtraction and multiplication of numbers. For
example, the commutative laws $a+b=b+a$ and $ab=ba$ are valid for polynomials
just as they are for numbers; the same holds for the associative laws $\left(
a+b\right)  +c=a+\left(  b+c\right)  $ and $\left(  ab\right)  c=a\left(
bc\right)  $ and the distributive laws $\left(  a+b\right)  c=ac+bc$ and
$a\left(  b+c\right)  =ab+ac$. Moreover, each polynomial $a$ satisfies
$a+\mathbf{0}=\mathbf{0}+a=a$ and $a\cdot\mathbf{0}=\mathbf{0}\cdot
a=\mathbf{0}$ and $a\cdot\mathbf{1}=\mathbf{1}\cdot a=a$ and $a+\left(
-a\right)  =\left(  -a\right)  +a=\mathbf{0}$.

Using the notations of Definition \ref{def.commring}, we can summarize this as
follows: The set $\mathbb{Q}\left[  X\right]  $, endowed with the operations
$+$ and $\cdot$ just defined, and with the elements $\mathbf{0}$ and
$\mathbf{1}$, is a commutative ring. It is called the \textit{(univariate)
polynomial ring over }$\mathbb{Q}$.

\textbf{(g)} Let $a=\left(  a_{0},a_{1},a_{2},\ldots\right)  \in
\mathbb{Q}\left[  X\right]  $ and $\lambda\in\mathbb{Q}$. Then, $\lambda a$
denotes the polynomial $\left(  \lambda a_{0},\lambda a_{1},\lambda
a_{2},\ldots\right)  \in\mathbb{Q}\left[  X\right]  $. (This equals the
polynomial $\left(  \operatorname*{const}\lambda\right)  \cdot a$; thus,
identifying $\lambda$ with $\operatorname*{const}\lambda$ does not cause any
inconsistencies here.)

\textbf{(h)} If $p=\left(  p_{0},p_{1},p_{2},\ldots\right)  \in\mathbb{Q}%
\left[  X\right]  $ is a nonzero polynomial, then the \textit{degree} of $p$
is defined to be the maximum $i\in\mathbb{N}$ satisfying $p_{i}\neq0$. If
$p\in\mathbb{Q}\left[  X\right]  $ is the zero polynomial, then the degree of
$p$ is defined to be $-\infty$. (Here, $-\infty$ is just a fancy symbol, not a
number.) The degree of a polynomial $p \in\mathbb{Q}\left[  X\right]  $ is
denoted $\deg p$. For example, $\deg\left(  0,4,0,-1,0,0,0,\ldots\right)  =3$.

\textbf{(i)} If $a=\left(  a_{0},a_{1},a_{2},\ldots\right)  \in\mathbb{Q}%
\left[  X\right]  $ and $n\in\mathbb{N}$, then a polynomial $a^{n}%
\in\mathbb{Q}\left[  X\right]  $ is defined to be the product
$\underbrace{aa\cdots a}_{n\text{ times}}$. (This is understood to be
$\mathbf{1}$ when $n=0$. In general, an empty product of polynomials is always
understood to be $\mathbf{1}$.)

\textbf{(j)} We let $X$ denote the polynomial $\left(  0,1,0,0,0,\ldots
\right)  \in\mathbb{Q}\left[  X\right]  $. (This is the polynomial whose
$1$-st coefficient is $1$ and whose other coefficients are $0$.) This
polynomial is called the \textit{indeterminate} of $\mathbb{Q}\left[
X\right]  $. It is easy to see that, for any $n\in\mathbb{N}$, we have%
\[
X^{n}=\left(  \underbrace{0,0,\ldots,0}_{n\text{ zeroes}},1,0,0,0,\ldots
\right)  .
\]

This polynomial $X$ finally provides an answer to the questions
\textquotedblleft what is an indeterminate\textquotedblright\ and
\textquotedblleft what is a formal expression\textquotedblright. Namely, let
$\left(  p_{0},p_{1},p_{2},\ldots\right)  \in\mathbb{Q}\left[  X\right]  $ be
any polynomial. Then, the sum $p_{0}+p_{1}X+p_{2}X^{2}+\cdots$ is well-defined
(it is an infinite sum, but due to (\ref{eq.def.polynomial-univar.finite}) it
has only finitely many nonzero addends), and it is easy to see that this sum
equals $\left(  p_{0},p_{1},p_{2},\ldots\right)  $. Thus,
\[
\left(  p_{0},p_{1},p_{2},\ldots\right)  =p_{0}+p_{1}X+p_{2}X^{2}%
+\cdots\ \ \ \ \ \ \ \ \ \ \text{for every }\left(  p_{0},p_{1},p_{2}%
,\ldots\right)  \in\mathbb{Q}\left[  X\right]  .
\]
This finally allows us to write a polynomial $\left(  p_{0},p_{1},p_{2}%
,\ldots\right)  $ as a sum $p_{0}+p_{1}X+p_{2}X^{2}+\cdots$ while remaining
honest; the sum $p_{0}+p_{1}X+p_{2}X^{2}+\cdots$ is no longer a
\textquotedblleft formal expression\textquotedblright\ of unclear meaning, nor
a function, but it is just an alternative way to write the sequence $\left(
p_{0},p_{1},p_{2},\ldots\right)  $. So, at last, our notion of a polynomial
resembles the intuitive notion of a polynomial!

Of course, we can write polynomials as finite sums as well. Indeed, if
$\left(  p_{0},p_{1},p_{2},\ldots\right)  \in\mathbb{Q}\left[  X\right]  $ is
a polynomial and $N$ is a nonnegative integer such that every $n>N$ satisfies
$p_{n}=0$, then%
\[
\left(  p_{0},p_{1},p_{2},\ldots\right)  =p_{0}+p_{1}X+p_{2}X^{2}+\cdots
=p_{0}+p_{1}X+\cdots+p_{N}X^{N}%
\]
(because addends can be discarded when they are $0$). For example,
\begin{align*}
\left(  4,1,0,0,0,\ldots\right)   &  =4+1X=4+X\ \ \ \ \ \ \ \ \ \ \text{and}\\
\left(  \dfrac{1}{2},0,\dfrac{1}{3},0,0,0,\ldots\right)   &  =\dfrac{1}%
{2}+0X+\dfrac{1}{3}X^{2}=\dfrac{1}{2}+\dfrac{1}{3}X^{2}.
\end{align*}

\textbf{(k)} For our definition of polynomials to be fully compatible with our
intuition, we are missing only one more thing: a way to evaluate a polynomial
at a number, or some other object (e.g., another polynomial or a function).
This is easy: Let $p=\left(  p_{0},p_{1},p_{2},\ldots\right)  \in
\mathbb{Q}\left[  X\right]  $ be a polynomial, and let $\alpha\in\mathbb{Q}$.
Then, $p\left(  \alpha\right)  $ means the number $p_{0}+p_{1}\alpha
+p_{2}\alpha^{2}+\cdots\in\mathbb{Q}$. (Again, the infinite sum $p_{0}%
+p_{1}\alpha+p_{2}\alpha^{2}+\cdots$ makes sense because of
(\ref{eq.def.polynomial-univar.finite}).) Similarly, we can define $p\left(
\alpha\right)  $ when $\alpha\in\mathbb{R}$ (but in this case, $p\left(
\alpha\right)  $ will be an element of $\mathbb{R}$) or when $\alpha
\in\mathbb{C}$ (in this case, $p\left(  \alpha\right)  \in\mathbb{C}$) or when
$\alpha$ is a square matrix with rational entries (in this case, $p\left(
\alpha\right)  $ will also be such a matrix) or when $\alpha$ is another
polynomial (in this case, $p\left(  \alpha\right)  $ is such a polynomial as well).

For example, if $p=\left(  1,-2,0,3,0,0,0,\ldots\right)  =1-2X+3X^{3}$, then
$p\left(  \alpha\right)  =1-2\alpha+3\alpha^{3}$ for every $\alpha$.

The map $\mathbb{Q}\rightarrow\mathbb{Q},\ \alpha\mapsto p\left(
\alpha\right)  $ is called the \textit{polynomial function described by }$p$.
As we said above, this function is not $p$, and it is not a good idea to
equate it with $p$.

If $\alpha$ is a number (or a square matrix, or another polynomial), then
$p\left(  \alpha\right)  $ is called the result of \textit{evaluating }$p$
\textit{at }$X=\alpha$ (or, simply, evaluating $p$ at $\alpha$), or the result
of \textit{substituting }$\alpha$\textit{ for }$X$\textit{ in }$p$. This
notation, of course, reminds of functions; nevertheless, (as we already said a
few times) $p$ is \textbf{not a function}.

Probably the simplest three cases of evaluation are the following ones:

\begin{itemize}
\item We have $p\left(  0\right)  =p_{0}+p_{1}0^{1}+p_{2}0^{2}+\cdots=p_{0}$.
In other words, evaluating $p$ at $X=0$ yields the constant term of $p$.

\item We have $p\left(  1\right)  =p_{0}+p_{1}1^{1}+p_{2}1^{2}+\cdots
=p_{0}+p_{1}+p_{2}+\cdots$. In other words, evaluating $p$ at $X=1$ yields the
sum of all coefficients of $p$.

\item We have $p\left(  X\right)  =p_{0}+p_{1}X^{1}+p_{2}X^{2}+\cdots
=p_{0}+p_{1}X+p_{2}X^{2}+\cdots=p$. In other words, evaluating $p$ at $X=X$
yields $p$ itself. This allows us to write $p\left(  X\right)  $ for $p$. Many
authors do so, just in order to stress that $p$ is a polynomial and that the
indeterminate is called $X$. It should be kept in mind that $X$ is \textbf{not
a variable} (just as $p$ is \textbf{not a function}); it is the (fixed!)
sequence $\left(  0,1,0,0,0,\ldots\right)  \in\mathbb{Q}\left[  X\right]  $
which serves as the indeterminate for polynomials in $\mathbb{Q}\left[
X\right]  $.
\end{itemize}

\textbf{(l)} Often, one wants (or is required) to give an indeterminate a name
other than $X$. (For instance, instead of polynomials with rational
coefficients, we could be considering polynomials whose coefficients
themselves are polynomials in $\mathbb{Q}\left[  X\right]  $; and then, we
would not be allowed to use the letter $X$ for the \textquotedblleft
new\textquotedblright\ indeterminate anymore, as it already means the
indeterminate of $\mathbb{Q}\left[  X\right]  $ !) This can be done, and the
rules are the following: Any letter (that does not already have a meaning) can
be used to denote the indeterminate; but then, the set of all polynomials has
to be renamed as $\mathbb{Q}\left[  \eta\right]  $, where $\eta$ is this
letter. For instance, if we want to denote the indeterminate as $x$, then we
have to denote the set by $\mathbb{Q}\left[  x\right]  $.

It is furthermore convenient to regard the sets $\mathbb{Q}\left[
\eta\right]  $ for different letters $\eta$ as distinct. Thus, for example,
the polynomial $3X^{2}+1$ is not the same as the polynomial $3Y^{2}+1$. (The
reason for doing so is that one sometimes wishes to view both of these
polynomials as polynomials in the two variables $X$ and $Y$.) Formally
speaking, this means that we should define a polynomial in $\mathbb{Q}\left[
\eta\right]  $ to be not just a sequence $\left(  p_{0},p_{1},p_{2}%
,\ldots\right)  $ of rational numbers, but actually a pair $\left(  \left(
p_{0},p_{1},p_{2},\ldots\right)  ,\text{\textquotedblleft}\eta
\text{\textquotedblright}\right)  $ of a sequence of rational numbers and the
letter $\eta$. (Here, \textquotedblleft$\eta$\textquotedblright\ really means
the letter $\eta$, not the sequence $\left(  0,1,0,0,0,\ldots\right)  $.) This
is, of course, a very technical point which is of little relevance to most of
mathematics; it becomes important when one tries to implement polynomials in a
programming language.

\textbf{(m)} As already explained, we can replace $\mathbb{Q}$ by $\mathbb{Z}%
$, $\mathbb{R}$, $\mathbb{C}$ or any other commutative ring $\mathbb{K}$ in
the above definition. (See Definition \ref{def.commring} for the definition of
a commutative ring.) When $\mathbb{Q}$ is replaced by a commutative ring
$\mathbb{K}$, the notion of \textquotedblleft univariate polynomials with
rational coefficients\textquotedblright\ becomes \textquotedblleft univariate
polynomials with coefficients in $\mathbb{K}$\textquotedblright\ (also known
as \textquotedblleft univariate polynomials over $\mathbb{K}$%
\textquotedblright), and the set of such polynomials is denoted by
$\mathbb{K}\left[  X\right]  $ rather than $\mathbb{Q}\left[  X\right]  $.
\end{definition}

So much for univariate polynomials.

Polynomials in multiple variables are (in my opinion) treated the best in
\cite[Chapter II, \S 3]{Lang02}, where they are introduced as elements of a
monoid ring. However, this treatment is rather abstract and uses a good deal
of algebraic language\footnote{Also, the book \cite{Lang02} is notorious for
its unpolished writing; it is best read with Bergman's companion
\cite{Bergman-Lang} at hand.}. The treatments in \cite[\S 4.5]{Walker87}, in
\cite[Chapter A-3]{Rotman15} and in \cite[Chapter IV, \S 4]{BirkMac} use the
above-mentioned recursive shortcut that makes them inferior (in my opinion). A
neat (and rather elementary) treatment of polynomials in $n$ variables (for
finite $n$) can be found in \cite[Chapter III, \S 5]{Hungerford-03}, in
\cite[\S 7.16]{Loehr-BC}, in \cite[\S 30.2]{GalQua18}, in \cite[\S 18]%
{ZarSam67} and in \cite[\S I.8]{AmaEsc05}; it generalizes the viewpoint we
used in Definition \ref{def.polynomial-univar} for univariate polynomials
above\footnote{You are reading right: The analysis textbook \cite{AmaEsc05} is
one of the few sources I am aware of to define the (algebraic!) notion of
polynomials precisely and well.}.

\section{\label{chp.ind}A closer look at induction}

In this chapter, we shall recall several versions of the \textit{induction
principle} (the principle of mathematical induction) and provide examples for
their use. We assume that the reader is at least somewhat familiar with
mathematical induction\footnote{If not, introductions can be found in
\cite[Chapter 5]{LeLeMe16}, \cite{Day-proofs}, \cite[Chapter 6]{Vellem06},
\cite[Chapter 10]{Hammac15}, \cite{Vorobi02} and various other sources.}; we
shall present some nonstandard examples of its use (including a proof of the
legitimacy of the definition of a sum $\sum_{s\in S}a_{s}$ given in Section
\ref{sect.sums-repetitorium}).

\subsection{\label{sect.ind.IP0}Standard induction}

\subsubsection{The Principle of Mathematical Induction}

We first recall the classical principle of mathematical
induction\footnote{Keep in mind that $\mathbb{N}$ means the set $\left\{
0,1,2,\ldots\right\}  $ for us.}:

\begin{theorem}
\label{thm.ind.IP0}For each $n\in\mathbb{N}$, let $\mathcal{A}\left(
n\right)  $ be a logical statement.

Assume the following:

\begin{statement}
\textit{Assumption 1:} The statement $\mathcal{A}\left(  0\right)  $ holds.
\end{statement}

\begin{statement}
\textit{Assumption 2:} If $m\in\mathbb{N}$ is such that $\mathcal{A}\left(
m\right)  $ holds, then $\mathcal{A}\left(  m+1\right)  $ also holds.
\end{statement}

Then, $\mathcal{A}\left(  n\right)  $ holds for each $n\in\mathbb{N}$.
\end{theorem}

Theorem \ref{thm.ind.IP0} is commonly taken to be one of the axioms of
mathematics (the \textquotedblleft axiom of induction\textquotedblright), or
(in type theory) as part of the definition of $\mathbb{N}$. Intuitively,
Theorem \ref{thm.ind.IP0} should be obvious: For example, if you want to prove
(under the assumptions of Theorem \ref{thm.ind.IP0}) that $\mathcal{A}\left(
4\right)  $ holds, you can argue as follows:

\begin{itemize}
\item By Assumption 1, the statement $\mathcal{A}\left(  0\right)  $ holds.

\item Thus, by Assumption 2 (applied to $m=0$), the statement $\mathcal{A}%
\left(  1\right)  $ holds.

\item Thus, by Assumption 2 (applied to $m=1$), the statement $\mathcal{A}%
\left(  2\right)  $ holds.

\item Thus, by Assumption 2 (applied to $m=2$), the statement $\mathcal{A}%
\left(  3\right)  $ holds.

\item Thus, by Assumption 2 (applied to $m=3$), the statement $\mathcal{A}%
\left(  4\right)  $ holds.
\end{itemize}

A similar (but longer) argument shows that the statement $\mathcal{A}\left(
5\right)  $ holds. Likewise, you can show that the statement $\mathcal{A}%
\left(  15\right)  $ holds, if you have the patience to apply Assumption 2 a
total of $15$ times. It is thus not surprising that $\mathcal{A}\left(
n\right)  $ holds for each $n\in\mathbb{N}$; but if you don't assume Theorem
\ref{thm.ind.IP0} as an axiom, you would need to write down a different proof
for each value of $n$ (which becomes the longer the larger $n$ is), and thus
would never reach the general result (i.e., that $\mathcal{A}\left(  n\right)
$ holds for \textbf{each} $n\in\mathbb{N}$), because you cannot write down
infinitely many proofs. What Theorem \ref{thm.ind.IP0} does is, roughly
speaking, to apply Assumption 2 for you as many times as it is needed for each
$n\in\mathbb{N}$.

(Authors of textbooks like to visualize Theorem \ref{thm.ind.IP0} by
envisioning an infinite sequence of dominos (numbered $0,1,2,\ldots$) placed
in row, sufficiently close to each other that if domino $m$ falls, then domino
$m+1$ will also fall. Now, assume that you kick domino $0$ over. What Theorem
\ref{thm.ind.IP0} then says is that each domino will fall. See, e.g.,
\cite[Chapter 10]{Hammac15} for a detailed explanation of this metaphor. Here
is another metaphor for Theorem \ref{thm.ind.IP0}: Assume that there is a
virus that infects nonnegative integers. Once it has infected some
$m\in\mathbb{N}$, it will soon spread to $m+1$ as well. Now, assume that $0$
gets infected. Then, Theorem \ref{thm.ind.IP0} says that each $n\in\mathbb{N}$
will eventually be infected.)

Theorem \ref{thm.ind.IP0} is called the \textit{principle of induction} or
\textit{principle of complete induction} or \textit{principle of mathematical
induction}, and we shall also call it \textit{principle of standard induction}
in order to distinguish it from several variant \textquotedblleft principles
of induction\textquotedblright\ that we will see later. Proofs that use this
principle are called \textit{proofs by induction} or \textit{induction
proofs}. Usually, in such proofs, we don't explicitly cite Theorem
\ref{thm.ind.IP0}, but instead say certain words that signal that Theorem
\ref{thm.ind.IP0} is being applied and that (ideally) also indicate what
statements $\mathcal{A}\left(  n\right)  $ it is being applied to\footnote{We
will explain this in Convention \ref{conv.ind.IP0lang} below.}. However, for
our very first example of a proof by induction, we are going to use Theorem
\ref{thm.ind.IP0} explicitly. We shall show the following fact:

\begin{proposition}
\label{prop.ind.ari-geo}Let $q$ and $d$ be two real numbers such that $q\neq
1$. Let $\left(  a_{0},a_{1},a_{2},\ldots\right)  $ be a sequence of real
numbers. Assume that%
\begin{equation}
a_{n+1}=qa_{n}+d\ \ \ \ \ \ \ \ \ \ \text{for each }n\in\mathbb{N}.
\label{eq.prop.ind.ari-geo.ass}%
\end{equation}
Then,
\begin{equation}
a_{n}=q^{n}a_{0}+\dfrac{q^{n}-1}{q-1}d\ \ \ \ \ \ \ \ \ \ \text{for each }%
n\in\mathbb{N}. \label{eq.prop.ind.ari-geo.claim}%
\end{equation}

\end{proposition}

\begin{proof}
[Proof of Proposition \ref{prop.ind.ari-geo}.]For each $n\in\mathbb{N}$, we
let $\mathcal{A}\left(  n\right)  $ be the statement \newline$\left(
a_{n}=q^{n}a_{0}+\dfrac{q^{n}-1}{q-1}d\right)  $. Thus, our goal is to prove
the statement $\mathcal{A}\left(  n\right)  $ for each $n\in\mathbb{N}$.

We first notice that the statement $\mathcal{A}\left(  0\right)  $
holds\footnote{\textit{Proof.} This is easy to verify: We have $q^{0}=1$, thus
$q^{0}-1=0$, and therefore $\dfrac{q^{0}-1}{q-1}=\dfrac{0}{q-1}=0$. Now,%
\[
\underbrace{q^{0}}_{=1}a_{0}+\underbrace{\dfrac{q^{0}-1}{q-1}}_{=0}%
d=1a_{0}+0d=a_{0},
\]
so that $a_{0}=q^{0}a_{0}+\dfrac{q^{0}-1}{q-1}d$. But this is precisely the
statement $\mathcal{A}\left(  0\right)  $ (since $\mathcal{A}\left(  0\right)
$ is defined to be the statement $\left(  a_{0}=q^{0}a_{0}+\dfrac{q^{0}%
-1}{q-1}d\right)  $). Hence, the statement $\mathcal{A}\left(  0\right)  $
holds.}.

Now, we claim that
\begin{equation}
\text{if }m\in\mathbb{N}\text{ is such that }\mathcal{A}\left(  m\right)
\text{ holds, then }\mathcal{A}\left(  m+1\right)  \text{ also holds.}
\label{pf.prop.ind.ari-geo.step}%
\end{equation}

[\textit{Proof of (\ref{pf.prop.ind.ari-geo.step}):} Let $m\in\mathbb{N}$ be
such that $\mathcal{A}\left(  m\right)  $ holds. We must show that
$\mathcal{A}\left(  m+1\right)  $ also holds.

We have assumed that $\mathcal{A}\left(  m\right)  $ holds. In other words,
$a_{m}=q^{m}a_{0}+\dfrac{q^{m}-1}{q-1}d$ holds\footnote{because $\mathcal{A}%
\left(  m\right)  $ is defined to be the statement $\left(  a_{m}=q^{m}%
a_{0}+\dfrac{q^{m}-1}{q-1}d\right)  $}. Now, (\ref{eq.prop.ind.ari-geo.ass})
(applied to $n=m$) yields%
\begin{align*}
a_{m+1}  &  =q\underbrace{a_{m}}_{=q^{m}a_{0}+\dfrac{q^{m}-1}{q-1}%
d}+d=q\left(  q^{m}a_{0}+\dfrac{q^{m}-1}{q-1}d\right)  +d\\
&  =\underbrace{qq^{m}}_{=q^{m+1}}a_{0}+\underbrace{q\cdot\dfrac{q^{m}-1}%
{q-1}d+d}_{=\left(  q\cdot\dfrac{q^{m}-1}{q-1}+1\right)  d}\\
&  =q^{m+1}a_{0}+\underbrace{\left(  q\cdot\dfrac{q^{m}-1}{q-1}+1\right)
}_{\substack{=\dfrac{q\left(  q^{m}-1\right)  +\left(  q-1\right)  }%
{q-1}=\dfrac{q^{m+1}-1}{q-1}\\\text{(since }q\left(  q^{m}-1\right)  +\left(
q-1\right)  =qq^{m}-q+q-1=qq^{m}-1=q^{m+1}-1\text{)}}}d\\
&  =q^{m+1}a_{0}+\dfrac{q^{m+1}-1}{q-1}d.
\end{align*}
So we have shown that $a_{m+1}=q^{m+1}a_{0}+\dfrac{q^{m+1}-1}{q-1}d$. But this
is precisely the statement $\mathcal{A}\left(  m+1\right)  $%
\ \ \ \ \footnote{because $\mathcal{A}\left(  m+1\right)  $ is defined to be
the statement $\left(  a_{m+1}=q^{m+1}a_{0}+\dfrac{q^{m+1}-1}{q-1}d\right)  $%
}. Thus, the statement $\mathcal{A}\left(  m+1\right)  $ holds.

Now, forget that we fixed $m$. We thus have shown that if $m\in\mathbb{N}$ is
such that $\mathcal{A}\left(  m\right)  $ holds, then $\mathcal{A}\left(
m+1\right)  $ also holds. This proves (\ref{pf.prop.ind.ari-geo.step}).]

Now, both assumptions of Theorem \ref{thm.ind.IP0} are satisfied (indeed,
Assumption 1 holds because the statement $\mathcal{A}\left(  0\right)  $
holds, whereas Assumption 2 holds because of (\ref{pf.prop.ind.ari-geo.step}%
)). Thus, Theorem \ref{thm.ind.IP0} shows that $\mathcal{A}\left(  n\right)  $
holds for each $n\in\mathbb{N}$. In other words, $a_{n}=q^{n}a_{0}%
+\dfrac{q^{n}-1}{q-1}d$ holds for each $n\in\mathbb{N}$ (since $\mathcal{A}%
\left(  n\right)  $ is the statement $\left(  a_{n}=q^{n}a_{0}+\dfrac{q^{n}%
-1}{q-1}d\right)  $). This proves Proposition \ref{prop.ind.ari-geo}.
\end{proof}

\subsubsection{Conventions for writing induction proofs}

Now, let us introduce some standard language that is commonly used in proofs
by induction:

\begin{convention}
\label{conv.ind.IP0lang}For each $n\in\mathbb{N}$, let $\mathcal{A}\left(
n\right)  $ be a logical statement. Assume that you want to prove that
$\mathcal{A}\left(  n\right)  $ holds for each $n\in\mathbb{N}$.

Theorem \ref{thm.ind.IP0} offers the following strategy for proving this:
First show that Assumption 1 of Theorem \ref{thm.ind.IP0} is satisfied; then,
show that Assumption 2 of Theorem \ref{thm.ind.IP0} is satisfied; then,
Theorem \ref{thm.ind.IP0} automatically completes your proof.

A proof that follows this strategy is called a \textit{proof by induction on
}$n$ (or \textit{proof by induction over }$n$) or (less precisely) an
\textit{inductive proof}. When you follow this strategy, you say that you are
\textit{inducting on }$n$ (or \textit{over }$n$). The proof that Assumption 1
is satisfied is called the \textit{induction base} (or \textit{base case}) of
the proof. The proof that Assumption 2 is satisfied is called the
\textit{induction step} of the proof.

In order to prove that Assumption 2 is satisfied, you will usually want to fix
an $m\in\mathbb{N}$ such that $\mathcal{A}\left(  m\right)  $ holds, and then
prove that $\mathcal{A}\left(  m+1\right)  $ holds. In other words, you will
usually want to fix $m\in\mathbb{N}$, assume that $\mathcal{A}\left(
m\right)  $ holds, and then prove that $\mathcal{A}\left(  m+1\right)  $
holds. When doing so, it is common to refer to the assumption that
$\mathcal{A}\left(  m\right)  $ holds as the \textit{induction hypothesis} (or
\textit{induction assumption}).
\end{convention}

Using this language, we can rewrite our above proof of Proposition
\ref{prop.ind.ari-geo} as follows:

\begin{proof}
[Proof of Proposition \ref{prop.ind.ari-geo} (second version).]For each
$n\in\mathbb{N}$, we let $\mathcal{A}\left(  n\right)  $ be the statement
$\left(  a_{n}=q^{n}a_{0}+\dfrac{q^{n}-1}{q-1}d\right)  $. Thus, our goal is
to prove the statement $\mathcal{A}\left(  n\right)  $ for each $n\in
\mathbb{N}$.

We shall prove this by induction on $n$:

\textit{Induction base:} We have $q^{0}=1$, thus $q^{0}-1=0$, and therefore
$\dfrac{q^{0}-1}{q-1}=\dfrac{0}{q-1}=0$. Now,%
\[
\underbrace{q^{0}}_{=1}a_{0}+\underbrace{\dfrac{q^{0}-1}{q-1}}_{=0}%
d=1a_{0}+0d=a_{0},
\]
so that $a_{0}=q^{0}a_{0}+\dfrac{q^{0}-1}{q-1}d$. But this is precisely the
statement $\mathcal{A}\left(  0\right)  $ (since $\mathcal{A}\left(  0\right)
$ is defined to be the statement $\left(  a_{0}=q^{0}a_{0}+\dfrac{q^{0}%
-1}{q-1}d\right)  $). Hence, the statement $\mathcal{A}\left(  0\right)  $
holds. This completes the induction base.

\textit{Induction step:} Let $m\in\mathbb{N}$. Assume that $\mathcal{A}\left(
m\right)  $ holds. We must show that $\mathcal{A}\left(  m+1\right)  $ also holds.

We have assumed that $\mathcal{A}\left(  m\right)  $ holds (this is our
induction hypothesis). In other words, $a_{m}=q^{m}a_{0}+\dfrac{q^{m}-1}%
{q-1}d$ holds\footnote{because $\mathcal{A}\left(  m\right)  $ is defined to
be the statement $\left(  a_{m}=q^{m}a_{0}+\dfrac{q^{m}-1}{q-1}d\right)  $}.
Now, (\ref{eq.prop.ind.ari-geo.ass}) (applied to $n=m$) yields%
\begin{align*}
a_{m+1}  &  =q\underbrace{a_{m}}_{=q^{m}a_{0}+\dfrac{q^{m}-1}{q-1}%
d}+d=q\left(  q^{m}a_{0}+\dfrac{q^{m}-1}{q-1}d\right)  +d\\
&  =\underbrace{qq^{m}}_{=q^{m+1}}a_{0}+\underbrace{q\cdot\dfrac{q^{m}-1}%
{q-1}d+d}_{=\left(  q\cdot\dfrac{q^{m}-1}{q-1}+1\right)  d}\\
&  =q^{m+1}a_{0}+\underbrace{\left(  q\cdot\dfrac{q^{m}-1}{q-1}+1\right)
}_{\substack{=\dfrac{q\left(  q^{m}-1\right)  +\left(  q-1\right)  }%
{q-1}=\dfrac{q^{m+1}-1}{q-1}\\\text{(since }q\left(  q^{m}-1\right)  +\left(
q-1\right)  =qq^{m}-q+q-1=qq^{m}-1=q^{m+1}-1\text{)}}}d\\
&  =q^{m+1}a_{0}+\dfrac{q^{m+1}-1}{q-1}d.
\end{align*}
So we have shown that $a_{m+1}=q^{m+1}a_{0}+\dfrac{q^{m+1}-1}{q-1}d$. But this
is precisely the statement $\mathcal{A}\left(  m+1\right)  $%
\ \ \ \ \footnote{because $\mathcal{A}\left(  m+1\right)  $ is defined to be
the statement $\left(  a_{m+1}=q^{m+1}a_{0}+\dfrac{q^{m+1}-1}{q-1}d\right)  $%
}. Thus, the statement $\mathcal{A}\left(  m+1\right)  $ holds.

Now, forget that we fixed $m$. We thus have shown that if $m\in\mathbb{N}$ is
such that $\mathcal{A}\left(  m\right)  $ holds, then $\mathcal{A}\left(
m+1\right)  $ also holds. This completes the induction step.

Thus, we have completed both the induction base and the induction step. Hence,
by induction, we conclude that $\mathcal{A}\left(  n\right)  $ holds for each
$n\in\mathbb{N}$. This proves Proposition \ref{prop.ind.ari-geo}.
\end{proof}

The proof we just gave still has a lot of \textquotedblleft
boilerplate\textquotedblright\ text. For example, we have explicitly defined
the statement $\mathcal{A}\left(  n\right)  $, but it is not really necessary,
since it is clear what this statement should be (viz., it should be the claim
we are proving, without the \textquotedblleft for each $n\in\mathbb{N}%
$\textquotedblright\ part). Allowing ourselves some imprecision, we could say
this statement is simply (\ref{eq.prop.ind.ari-geo.claim}). (This is a bit
imprecise, because (\ref{eq.prop.ind.ari-geo.claim}) contains the words
\textquotedblleft for each $n\in\mathbb{N}$\textquotedblright, but it should
be clear that we don't mean to include these words, since there can be no
\textquotedblleft for each $n\in\mathbb{N}$\textquotedblright\ in the
statement $\mathcal{A}\left(  n\right)  $.) Furthermore, we don't need to
write the sentence

\begin{quote}
\textquotedblleft Thus, we have completed both the induction base and the
induction step\textquotedblright
\end{quote}

\noindent before we declare our inductive proof to be finished; it is clear
enough that we have completed them. We also can remove the following two sentences:

\begin{quote}
\textquotedblleft Now, forget that we fixed $m$. We thus have shown that if
$m\in\mathbb{N}$ is such that $\mathcal{A}\left(  m\right)  $ holds, then
$\mathcal{A}\left(  m+1\right)  $ also holds.\textquotedblright.
\end{quote}

\noindent In fact, these sentences merely say that we have completed the
induction step; they carry no other information (since the induction step
always consists in fixing $m\in\mathbb{N}$ such that $\mathcal{A}\left(
m\right)  $ holds, and proving that $\mathcal{A}\left(  m+1\right)  $ also
holds). So once we say that the induction step is completed, we don't need
these sentences anymore.

So we can shorten our proof above a bit further:

\begin{proof}
[Proof of Proposition \ref{prop.ind.ari-geo} (third version).]We shall prove
(\ref{eq.prop.ind.ari-geo.claim}) by induction on $n$:

\textit{Induction base:} We have $q^{0}=1$, thus $q^{0}-1=0$, and therefore
$\dfrac{q^{0}-1}{q-1}=\dfrac{0}{q-1}=0$. Now,%
\[
\underbrace{q^{0}}_{=1}a_{0}+\underbrace{\dfrac{q^{0}-1}{q-1}}_{=0}%
d=1a_{0}+0d=a_{0},
\]
so that $a_{0}=q^{0}a_{0}+\dfrac{q^{0}-1}{q-1}d$. In other words,
(\ref{eq.prop.ind.ari-geo.claim}) holds for $n=0$.\ \ \ \ \footnote{Note that
the statement \textquotedblleft(\ref{eq.prop.ind.ari-geo.claim}) holds for
$n=0$\textquotedblright\ (which we just proved) is precisely the statement
$\mathcal{A}\left(  0\right)  $ in the previous two versions of our proof.}
This completes the induction base.

\textit{Induction step:} Let $m\in\mathbb{N}$. Assume that
(\ref{eq.prop.ind.ari-geo.claim}) holds for $n=m$.\ \ \ \ \footnote{Note that
the statement \textquotedblleft(\ref{eq.prop.ind.ari-geo.claim}) holds for
$n=m$\textquotedblright\ (which we just assumed) is precisely the statement
$\mathcal{A}\left(  m\right)  $ in the previous two versions of our proof.} We
must show that (\ref{eq.prop.ind.ari-geo.claim}) holds for $n=m+1$%
.\ \ \ \ \footnote{Note that this statement \textquotedblleft%
(\ref{eq.prop.ind.ari-geo.claim}) holds for $n=m+1$\textquotedblright\ is
precisely the statement $\mathcal{A}\left(  m+1\right)  $ in the previous two
versions of our proof.}

We have assumed that (\ref{eq.prop.ind.ari-geo.claim}) holds for $n=m$. In
other words, $a_{m}=q^{m}a_{0}+\dfrac{q^{m}-1}{q-1}d$ holds. Now,
(\ref{eq.prop.ind.ari-geo.ass}) (applied to $n=m$) yields%
\begin{align*}
a_{m+1}  &  =q\underbrace{a_{m}}_{=q^{m}a_{0}+\dfrac{q^{m}-1}{q-1}%
d}+d=q\left(  q^{m}a_{0}+\dfrac{q^{m}-1}{q-1}d\right)  +d\\
&  =\underbrace{qq^{m}}_{=q^{m+1}}a_{0}+\underbrace{q\cdot\dfrac{q^{m}-1}%
{q-1}d+d}_{=\left(  q\cdot\dfrac{q^{m}-1}{q-1}+1\right)  d}\\
&  =q^{m+1}a_{0}+\underbrace{\left(  q\cdot\dfrac{q^{m}-1}{q-1}+1\right)
}_{\substack{=\dfrac{q\left(  q^{m}-1\right)  +\left(  q-1\right)  }%
{q-1}=\dfrac{q^{m+1}-1}{q-1}\\\text{(since }q\left(  q^{m}-1\right)  +\left(
q-1\right)  =qq^{m}-q+q-1=qq^{m}-1=q^{m+1}-1\text{)}}}d\\
&  =q^{m+1}a_{0}+\dfrac{q^{m+1}-1}{q-1}d.
\end{align*}
So we have shown that $a_{m+1}=q^{m+1}a_{0}+\dfrac{q^{m+1}-1}{q-1}d$. In other
words, (\ref{eq.prop.ind.ari-geo.claim}) holds for $n=m+1$. This completes the
induction step. Hence, (\ref{eq.prop.ind.ari-geo.claim}) is proven by
induction. This proves Proposition \ref{prop.ind.ari-geo}.
\end{proof}

\subsection{Examples from modular arithmetic}

\subsubsection{Divisibility of integers}

We shall soon give some more examples of inductive proofs, including some that
will include slightly new tactics. These examples come from the realm of
\textit{modular arithmetic}, which is the study of congruences modulo
integers. Before we come to these examples, we will introduce the definition
of such congruences. But first, let us recall the definition of divisibility:

\begin{definition}
\label{def.divisibility}Let $u$ and $v$ be two integers. Then, we say that $u$
\textit{divides} $v$ if and only if there exists an integer $w$ such that
$v=uw$. Instead of saying \textquotedblleft$u$ divides $v$\textquotedblright,
we can also say \textquotedblleft$v$ is \textit{divisible by }$u$%
\textquotedblright\ or \textquotedblleft$v$ is a \textit{multiple} of
$u$\textquotedblright\ or \textquotedblleft$u$ is a \textit{divisor} of
$v$\textquotedblright\ or \textquotedblleft$u\mid v$\textquotedblright.
\end{definition}

Thus, two integers $u$ and $v$ satisfy $u\mid v$ if and only if there is some
$w\in\mathbb{Z}$ such that $v=uw$. For example, $1\mid v$ holds for every
integer $v$ (since $v=1v$), whereas $0\mid v$ holds only for $v=0$ (since
$v=0w$ is equivalent to $v=0$). An integer $v$ satisfies $2\mid v$ if and only
if $v$ is even.

Definition \ref{def.divisibility} is fairly common in the modern literature
(e.g., it is used in \cite{Day-proofs}, \cite{LeLeMe16}, \cite{Mulhol16} and
\cite{Rotman15}), but there are also some books that define these notations
differently. For example, in \cite{GKP}, the notation \textquotedblleft$u$
divides $v$\textquotedblright\ is defined differently (it requires not only
the existence of an integer $w$ such that $v=uw$, but also that $u$ is
positive), whereas the notation \textquotedblleft$v$ is a multiple of
$u$\textquotedblright\ is defined as it is here (i.e., it just means that
there exists an integer $w$ such that $v=uw$); thus, these two notations are
not mutually interchangeable in \cite{GKP}.

Let us first prove some basic properties of divisibility:

\begin{proposition}
\label{prop.div.trans}Let $a$, $b$ and $c$ be three integers such that $a\mid
b$ and $b\mid c$. Then, $a\mid c$.
\end{proposition}

\begin{proof}
[Proof of Proposition \ref{prop.div.trans}.]We have $a\mid b$. In other words,
there exists an integer $w$ such that $b=aw$ (by the definition of
\textquotedblleft divides\textquotedblright). Consider this $w$, and denote it
by $k$. Thus, $k$ is an integer such that $b=ak$.

We have $b\mid c$. In other words, there exists an integer $w$ such that
$c=bw$ (by the definition of \textquotedblleft divides\textquotedblright).
Consider this $w$, and denote it by $j$. Thus, $j$ is an integer such that
$c=bj$.

Now, $c=\underbrace{b}_{=ak}j=akj$. Hence, there exists an integer $w$ such
that $c=aw$ (namely, $w=kj$). In other words, $a$ divides $c$ (by the
definition of \textquotedblleft divides\textquotedblright). In other words,
$a\mid c$. This proves Proposition \ref{prop.div.trans}.
\end{proof}

\begin{proposition}
\label{prop.div.acbc}Let $a$, $b$ and $c$ be three integers such that $a\mid
b$. Then, $ac\mid bc$.
\end{proposition}

\begin{proof}
[Proof of Proposition \ref{prop.div.acbc}.]We have $a\mid b$. In other words,
there exists an integer $w$ such that $b=aw$ (by the definition of
\textquotedblleft divides\textquotedblright). Consider this $w$, and denote it
by $k$. Thus, $k$ is an integer such that $b=ak$. Hence, $\underbrace{b}%
_{=ak}c=akc=ack$. Thus, there exists an integer $w$ such that $bc=acw$
(namely, $w=k$). In other words, $ac$ divides $bc$ (by the definition of
\textquotedblleft divides\textquotedblright). In other words, $ac\mid bc$.
This proves Proposition \ref{prop.div.acbc}.
\end{proof}

\begin{proposition}
\label{prop.div.ax+by}Let $a$, $b$, $g$, $x$ and $y$ be integers such that
$g=ax+by$. Let $d$ be an integer such that $d\mid a$ and $d\mid b$. Then,
$d\mid g$.
\end{proposition}

\begin{proof}
[Proof of Proposition \ref{prop.div.ax+by}.]We have $d\mid a$. In other words,
there exists an integer $w$ such that $a=dw$ (by the definition of
\textquotedblleft divides\textquotedblright). Consider this $w$, and denote it
by $p$. Thus, $p$ is an integer and satisfies $a=dp$.

\begin{vershort}
Similarly, there is an integer $q$ such that $b=dq$. Consider this $q$.
\end{vershort}

\begin{verlong}
We have $d\mid b$. In other words, there exists an integer $w$ such that
$b=dw$ (by the definition of \textquotedblleft divides\textquotedblright).
Consider this $w$, and denote it by $q$. Thus, $q$ is an integer and satisfies
$b=dq$.
\end{verlong}

Now, $g=\underbrace{a}_{=dp}x+\underbrace{b}_{=dq}y=dpx+dqy=d\left(
px+qy\right)  $. Hence, there exists an integer $w$ such that $g=dw$ (namely,
$w=px+qy$). In other words, $d\mid g$ (by the definition of \textquotedblleft
divides\textquotedblright). This proves Proposition \ref{prop.div.ax+by}.
\end{proof}

It is easy to characterize divisibility in terms of fractions:

\begin{proposition}
\label{prop.div.frac}Let $a$ and $b$ be two integers such that $a\neq0$. Then,
$a\mid b$ if and only if $b/a$ is an integer.
\end{proposition}

\begin{proof}
[Proof of Proposition \ref{prop.div.frac}.]We first claim the following
logical implication\footnote{A \textit{logical implication} (or, short,
\textit{implication}) is a logical statement of the form \textquotedblleft if
$\mathcal{A}$, then $\mathcal{B}$\textquotedblright\ (where $\mathcal{A}$ and
$\mathcal{B}$ are two statements).}:%
\begin{equation}
\left(  a\mid b\right)  \ \Longrightarrow\ \left(  b/a\text{ is an
integer}\right)  . \label{pf.prop.div.frac.1}%
\end{equation}

[\textit{Proof of (\ref{pf.prop.div.frac.1}):} Assume that $a\mid b$. In other
words, there exists an integer $w$ such that $b=aw$ (by the definition of
\textquotedblleft divides\textquotedblright). Consider this $w$. Now, dividing
the equality $b=aw$ by $a$, we obtain $b/a=w$ (since $a\neq0$). Hence, $b/a$
is an integer (since $w$ is an integer). This proves the implication
(\ref{pf.prop.div.frac.1}).]

Next, we claim the following logical implication:%
\begin{equation}
\left(  b/a\text{ is an integer}\right)  \ \Longrightarrow\ \left(  a\mid
b\right)  . \label{pf.prop.div.frac.2}%
\end{equation}

[\textit{Proof of (\ref{pf.prop.div.frac.2}):} Assume that $b/a$ is an
integer. Let $k$ denote this integer. Thus, $b/a=k$, so that $b=ak$. Hence,
there exists an integer $w$ such that $b=aw$ (namely, $w=k$). In other words,
$a$ divides $b$ (by the definition of \textquotedblleft
divides\textquotedblright). In other words, $a\mid b$. This proves the
implication (\ref{pf.prop.div.frac.2}).]

Combining the implications (\ref{pf.prop.div.frac.1}) and
(\ref{pf.prop.div.frac.2}), we obtain the equivalence $\left(  a\mid b\right)
\ \Longleftrightarrow\ \left(  b/a\text{ is an integer}\right)  $. In other
words, $a\mid b$ if and only if $b/a$ is an integer. This proves Proposition
\ref{prop.div.frac}.
\end{proof}

\subsubsection{Definition of congruences}

We can now define congruences:

\begin{definition}
\label{def.mod.equiv}Let $a$, $b$ and $n$ be three integers. Then, we say that
$a$\textit{ is congruent to }$b$\textit{ modulo }$n$ if and only if $n\mid
a-b$. We shall use the notation \textquotedblleft$a\equiv b\operatorname{mod}%
n$\textquotedblright\ for \textquotedblleft$a$ is congruent to $b$ modulo
$n$\textquotedblright. Relations of the form \textquotedblleft$a\equiv
b\operatorname{mod}n$\textquotedblright\ (for integers $a$, $b$ and $n$) are
called \textit{congruences modulo }$n$.
\end{definition}

Thus, three integers $a$, $b$ and $n$ satisfy $a\equiv b\operatorname{mod}n$
if and only if $n\mid a-b$.

Hence, in particular:

\begin{itemize}
\item Any two integers $a$ and $b$ satisfy $a\equiv b\operatorname{mod}1$.
(Indeed, any two integers $a$ and $b$ satisfy $a-b=1\left(  a-b\right)  $,
thus $1\mid a-b$, thus $a\equiv b\operatorname{mod}1$.)

\item Two integers $a$ and $b$ satisfy $a\equiv b\operatorname{mod}0$ if and
only if $a=b$. (Indeed, $a\equiv b\operatorname{mod}0$ is equivalent to $0\mid
a-b$, which in turn is equivalent to $a-b=0$, which in turn is equivalent to
$a=b$.)

\item Two integers $a$ and $b$ satisfy $a\equiv b\operatorname{mod}2$ if and
only if they have the same parity (i.e., they are either both odd or both
even). This is not obvious at this point yet, but follows from Proposition
\ref{prop.mod.parity} further below.
\end{itemize}

We have%
\[
4\equiv10\operatorname{mod}3\ \ \ \ \ \ \ \ \ \ \text{and}%
\ \ \ \ \ \ \ \ \ \ 5\equiv-35\operatorname{mod}4.
\]

Note that Day, in \cite{Day-proofs}, writes \textquotedblleft$a\equiv_{n}%
b$\textquotedblright\ instead of \textquotedblleft$a\equiv b\operatorname{mod}%
n$\textquotedblright. Also, other authors (particularly of older texts) write
\textquotedblleft$a\equiv b\pmod{n}$\textquotedblright\ instead of
\textquotedblleft$a\equiv b\operatorname{mod}n$\textquotedblright.

Let us next introduce notations for the negations of the statements
\textquotedblleft$u\mid v$\textquotedblright\ and \textquotedblleft$a\equiv
b\operatorname{mod}n$\textquotedblright:

\begin{definition}
\textbf{(a)} If $u$ and $v$ are two integers, then the notation
\textquotedblleft$u\nmid v$\textquotedblright\ shall mean \textquotedblleft
not $u\mid v$\textquotedblright\ (that is, \textquotedblleft$u$ does not
divide $v$\textquotedblright).

\textbf{(b)} If $a$, $b$ and $n$ are three integers, then the notation
\textquotedblleft$a\not \equiv b\operatorname{mod}n$\textquotedblright\ shall
mean \textquotedblleft not $a\equiv b\operatorname{mod}n$\textquotedblright%
\ (that is, \textquotedblleft$a$ is not congruent to $b$ modulo $n$%
\textquotedblright).
\end{definition}

Thus, three integers $a$, $b$ and $n$ satisfy $a\not \equiv
b\operatorname{mod}n$ if and only if $n\nmid a-b$. For example, $1\not \equiv
-1\operatorname{mod}3$, since $3\nmid1-\left(  -1\right)  $.

\subsubsection{Congruence basics}

Let us now state some of the basic laws of congruences (so far, not needing
induction to prove):

\begin{proposition}
\label{prop.mod.0}Let $a$ and $n$ be integers. Then:

\textbf{(a)} We have $a\equiv0\operatorname{mod}n$ if and only if $n\mid a$.

\textbf{(b)} Let $b$ be an integer. Then, $a\equiv b\operatorname{mod}n$ if
and only if $a\equiv b\operatorname{mod}\left(  -n\right)  $.

\textbf{(c)} Let $m$ and $b$ be integers such that $m\mid n$. If $a\equiv
b\operatorname{mod}n$, then $a\equiv b\operatorname{mod}m$.
\end{proposition}

\begin{proof}
[Proof of Proposition \ref{prop.mod.0}.]\textbf{(a)} We have the following
chain of logical equivalences:%
\begin{align*}
&  \ \left(  a\equiv0\operatorname{mod}n\right) \\
&  \Longleftrightarrow\ \left(  a\text{ is congruent to }0\text{ modulo
}n\right) \\
&  \ \ \ \ \ \ \ \ \ \ \left(  \text{since \textquotedblleft}a\equiv
0\operatorname{mod}n\text{\textquotedblright\ is just a notation for
\textquotedblleft}a\text{ is congruent to }0\text{ modulo }%
n\text{\textquotedblright}\right) \\
&  \Longleftrightarrow\ \left(  n\mid\underbrace{a-0}_{=a}\right)
\ \ \ \ \ \ \ \ \ \ \left(  \text{by the definition of \textquotedblleft
congruent\textquotedblright}\right) \\
&  \Longleftrightarrow\ \left(  n\mid a\right)  .
\end{align*}
Thus, we have $a\equiv0\operatorname{mod}n$ if and only if $n\mid a$. This
proves Proposition \ref{prop.mod.0} \textbf{(a)}.

\textbf{(b)} Let us first assume that $a\equiv b\operatorname{mod}n$. Thus,
$a$ is congruent to $b$ modulo $n$. In other words, $n\mid a-b$ (by the
definition of \textquotedblleft congruent\textquotedblright). In other words,
$n$ divides $a-b$. In other words, there exists an integer $w$ such that
$a-b=nw$ (by the definition of \textquotedblleft divides\textquotedblright).
Consider this $w$, and denote it by $k$. Thus, $k$ is an integer such that
$a-b=nk$.

Thus, $a-b=nk=\left(  -n\right)  \left(  -k\right)  $. Hence, there exists an
integer $w$ such that $a-b=\left(  -n\right)  w$ (namely, $w=-k$). In other
words, $-n$ divides $a-b$ (by the definition of \textquotedblleft
divides\textquotedblright). In other words, $-n\mid a-b$. In other words, $a$
is congruent to $b$ modulo $-n$ (by the definition of \textquotedblleft
congruent\textquotedblright). In other words, $a\equiv b\operatorname{mod}%
\left(  -n\right)  $.

Now, forget that we assumed that $a\equiv b\operatorname{mod}n$. We thus have
shown that%
\begin{equation}
\text{if }a\equiv b\operatorname{mod}n\text{, then }a\equiv
b\operatorname{mod}\left(  -n\right)  . \label{pf.prop.mod.0.b.1}%
\end{equation}
The same argument (applied to $-n$ instead of $n$) shows that%
\[
\text{if }a\equiv b\operatorname{mod}\left(  -n\right)  \text{, then }a\equiv
b\operatorname{mod}\left(  -\left(  -n\right)  \right)  .
\]
Since $-\left(  -n\right)  =n$, this rewrites as follows:%
\[
\text{if }a\equiv b\operatorname{mod}\left(  -n\right)  \text{, then }a\equiv
b\operatorname{mod}n.
\]
Combining this implication with (\ref{pf.prop.mod.0.b.1}), we conclude that
$a\equiv b\operatorname{mod}n$ if and only if $a\equiv b\operatorname{mod}%
\left(  -n\right)  $. This proves Proposition \ref{prop.mod.0} \textbf{(b)}.

\textbf{(c)} Assume that $a\equiv b\operatorname{mod}n$. Thus, $a$ is
congruent to $b$ modulo $n$. In other words, $n\mid a-b$ (by the definition of
\textquotedblleft congruent\textquotedblright). Hence, Proposition
\ref{prop.div.trans} (applied to $m$, $n$ and $a-b$ instead of $a$, $b$ and
$c$) yields $m\mid a-b$ (since $m\mid n$). In other words, $a$ is congruent to
$b$ modulo $m$ (by the definition of \textquotedblleft
congruent\textquotedblright). Thus, $a\equiv b\operatorname{mod}m$. This
proves Proposition \ref{prop.mod.0} \textbf{(c)}.
\end{proof}

\begin{proposition}
\label{prop.mod.transi}Let $n$ be an integer.

\textbf{(a)} For any integer $a$, we have $a\equiv a\operatorname{mod}n$.

\textbf{(b)} For any integers $a$ and $b$ satisfying $a\equiv
b\operatorname{mod}n$, we have $b\equiv a\operatorname{mod}n$.

\textbf{(c)} For any integers $a$, $b$ and $c$ satisfying $a\equiv
b\operatorname{mod}n$ and $b\equiv c\operatorname{mod}n$, we have $a\equiv
c\operatorname{mod}n$.
\end{proposition}

\begin{proof}
[Proof of Proposition \ref{prop.mod.transi}.]\textbf{(a)} Let $a$ be an
integer. Then, $a-a=0=n\cdot0$. Hence, there exists an integer $w$ such that
$a-a=nw$ (namely, $w=0$). In other words, $n$ divides $a-a$ (by the definition
of \textquotedblleft divides\textquotedblright). In other words, $n\mid a-a$.
In other words, $a$ is congruent to $a$ modulo $n$ (by the definition of
\textquotedblleft congruent\textquotedblright). In other words, $a\equiv
a\operatorname{mod}n$. This proves Proposition \ref{prop.mod.transi}
\textbf{(a)}.

\textbf{(b)} Let $a$ and $b$ be two integers satisfying $a\equiv
b\operatorname{mod}n$. Thus, $a$ is congruent to $b$ modulo $n$ (since
$a\equiv b\operatorname{mod}n$). In other words, $n\mid a-b$ (by the
definition of \textquotedblleft congruent\textquotedblright). In other words,
$n$ divides $a-b$. In other words, there exists an integer $w$ such that
$a-b=nw$ (by the definition of \textquotedblleft divides\textquotedblright).
Consider this $w$, and denote it by $q$. Thus, $q$ is an integer such that
$a-b=nq$. Now, $b-a=-\underbrace{\left(  a-b\right)  }_{=nq}=-nq=n\left(
-q\right)  $. Hence, there exists an integer $w$ such that $b-a=nw$ (namely,
$w=-q$). In other words, $n$ divides $b-a$ (by the definition of
\textquotedblleft divides\textquotedblright). In other words, $n\mid b-a$. In
other words, $b$ is congruent to $a$ modulo $n$ (by the definition of
\textquotedblleft congruent\textquotedblright). In other words, $b\equiv
a\operatorname{mod}n$. This proves Proposition \ref{prop.mod.transi}
\textbf{(b)}.

\textbf{(c)} Let $a$, $b$ and $c$ be three integers satisfying $a\equiv
b\operatorname{mod}n$ and $b\equiv c\operatorname{mod}n$.

\begin{vershort}
Just as in the above proof of Proposition \ref{prop.mod.transi} \textbf{(b)},
we can use the assumption $a\equiv b\operatorname{mod}n$ to construct an
integer $q$ such that $a-b=nq$. Similarly, we can use the assumption $b\equiv
c\operatorname{mod}n$ to construct an integer $r$ such that $b-c=nr$. Consider
these $q$ and $r$.
\end{vershort}

\begin{verlong}
From $a\equiv b\operatorname{mod}n$, we conclude that $a$ is congruent to $b$
modulo $n$. In other words, $n\mid a-b$ (by the definition of
\textquotedblleft congruent\textquotedblright). In other words, $n$ divides
$a-b$. In other words, there exists an integer $w$ such that $a-b=nw$ (by the
definition of \textquotedblleft divides\textquotedblright). Consider this $w$,
and denote it by $q$. Thus, $q$ is an integer such that $a-b=nq$.

From $b\equiv c\operatorname{mod}n$, we conclude that $b$ is congruent to $c$
modulo $n$. In other words, $n\mid b-c$ (by the definition of
\textquotedblleft congruent\textquotedblright). In other words, $n$ divides
$b-c$. In other words, there exists an integer $w$ such that $b-c=nw$ (by the
definition of \textquotedblleft divides\textquotedblright). Consider this $w$,
and denote it by $r$. Thus, $r$ is an integer such that $b-c=nr$.
\end{verlong}

Now,%
\[
a-c=\underbrace{\left(  a-b\right)  }_{=nq}+\underbrace{\left(  b-c\right)
}_{=nr}=nq+nr=n\left(  q+r\right)  .
\]
Hence, there exists an integer $w$ such that $a-c=nw$ (namely, $w=q+r$). In
other words, $n$ divides $a-c$ (by the definition of \textquotedblleft
divides\textquotedblright). In other words, $n\mid a-c$. In other words, $a$
is congruent to $c$ modulo $n$ (by the definition of \textquotedblleft
congruent\textquotedblright). In other words, $a\equiv c\operatorname{mod}n$.
This proves Proposition \ref{prop.mod.transi} \textbf{(c)}.
\end{proof}

Simple as they are, the three parts of Proposition \ref{prop.mod.transi} have
names: Proposition \ref{prop.mod.transi} \textbf{(a)} is called the
\textit{reflexivity of congruence (modulo }$n$\textit{)}; Proposition
\ref{prop.mod.transi} \textbf{(b)} is called the \textit{symmetry of
congruence (modulo }$n$\textit{)}; Proposition \ref{prop.mod.transi}
\textbf{(c)} is called the \textit{transitivity of congruence (modulo }%
$n$\textit{)}.

Proposition \ref{prop.mod.transi} \textbf{(b)} allows the following definition:

\begin{definition}
Let $n$, $a$ and $b$ be three integers. Then, we say that $a$ \textit{and }$b$
\textit{are congruent modulo }$n$ if and only if $a\equiv b\operatorname{mod}%
n$. Proposition \ref{prop.mod.transi} \textbf{(b)} shows that $a$ and $b$
actually play equal roles in this relation (i.e., the statement
\textquotedblleft$a$ and $b$ are congruent modulo $n$\textquotedblright\ is
equivalent to \textquotedblleft$b$ and $a$ are congruent modulo $n$%
\textquotedblright).
\end{definition}

\begin{proposition}
\label{prop.mod.n=0}Let $n$ be an integer. Then, $n\equiv0\operatorname{mod}n$.
\end{proposition}

\begin{proof}
[Proof of Proposition \ref{prop.mod.n=0}.]We have $n=n\cdot1$. Thus, there
exists an integer $w$ such that $n=nw$ (namely, $w=1$). Therefore, $n\mid n$
(by the definition of \textquotedblleft divides\textquotedblright).
Proposition \ref{prop.mod.0} \textbf{(a)} (applied to $a=n$) shows that we
have $n\equiv0\operatorname{mod}n$ if and only if $n\mid n$. Hence, we have
$n\equiv0\operatorname{mod}n$ (since $n\mid n$). This proves Proposition
\ref{prop.mod.n=0}.
\end{proof}

\subsubsection{Chains of congruences}

Proposition \ref{prop.mod.transi} shows that congruences (modulo $n$) behave
like equalities -- in that we can turn them around (since Proposition
\ref{prop.mod.transi} \textbf{(b)} says that $a\equiv b\operatorname{mod}n$
implies $b\equiv a\operatorname{mod}n$) and we can chain them together (by
Proposition \ref{prop.mod.transi} \textbf{(c)}) and in that every integer is
congruent to itself (by Proposition \ref{prop.mod.transi} \textbf{(a)}). This
leads to the following notation:

\begin{definition}
If $a_{1},a_{2},\ldots,a_{k}$ and $n$ are integers, then the statement
\textquotedblleft$a_{1}\equiv a_{2}\equiv\cdots\equiv a_{k}\operatorname{mod}%
n$\textquotedblright\ shall mean that
\[
\left(  a_{i}\equiv a_{i+1}\operatorname{mod}n\text{ holds for each }%
i\in\left\{  1,2,\ldots,k-1\right\}  \right)  .
\]
Such a statement is called a \textit{chain of congruences modulo }$n$ (or,
less precisely, a \textit{chain of congruences}). We shall refer to the
integers $a_{1},a_{2},\ldots,a_{k}$ (but not $n$) as the \textit{members} of
this chain.
\end{definition}

For example, the chain $a\equiv b\equiv c\equiv d\operatorname{mod}n$ (for
five integers $a,b,c,d,n$) means that $a\equiv b\operatorname{mod}n$ and
$b\equiv c\operatorname{mod}n$ and $c\equiv d\operatorname{mod}n$.

The usefulness of such chains lies in the following fact:

\begin{proposition}
\label{prop.mod.chain}Let $a_{1},a_{2},\ldots,a_{k}$ and $n$ be integers such
that $a_{1}\equiv a_{2}\equiv\cdots\equiv a_{k}\operatorname{mod}n$. Let $u$
and $v$ be two elements of $\left\{  1,2,\ldots,k\right\}  $. Then,%
\[
a_{u}\equiv a_{v}\operatorname{mod}n.
\]

\end{proposition}

In other words, any two members of a chain of congruences modulo $n$ are
congruent to each other modulo $n$. Thus, chains of congruences are like
chains of equalities: From any chain of congruences modulo $n$ with $k$
members, you can extract $k^{2}$ congruences modulo $n$ by picking any two
members of the chain.

\begin{example}
Proposition \ref{prop.mod.chain} shows (among other things) that if
$a,b,c,d,e,n$ are integers such that $a\equiv b\equiv c\equiv d\equiv
e\operatorname{mod}n$, then $a\equiv d\operatorname{mod}n$ and $b\equiv
d\operatorname{mod}n$ and $e\equiv b\operatorname{mod}n$ (and various other congruences).
\end{example}

Unsurprisingly, Proposition \ref{prop.mod.chain} can be proven by induction,
although not in an immediately obvious manner: We cannot directly prove it by
induction on $n$, on $k$, on $u$ or on $v$. Instead, we will first introduce
an auxiliary statement (the statement (\ref{pf.prop.mod.chain.Ai}) in the
following proof) which will be tailored to an inductive proof. This is a
commonly used tactic, and particularly helpful to us now as we only have the
most basic form of the principle of induction available. (Soon, we will see
more versions of that principle, which will obviate the need for some of the tailoring.)

\begin{proof}
[Proof of Proposition \ref{prop.mod.chain}.]By assumption, we have
$a_{1}\equiv a_{2}\equiv\cdots\equiv a_{k}\operatorname{mod}n$. In other
words,
\begin{equation}
\left(  a_{i}\equiv a_{i+1}\operatorname{mod}n\text{ holds for each }%
i\in\left\{  1,2,\ldots,k-1\right\}  \right)  \label{pf.prop.mod.chain.ass}%
\end{equation}
(since this is what \textquotedblleft$a_{1}\equiv a_{2}\equiv\cdots\equiv
a_{k}\operatorname{mod}n$\textquotedblright\ means).

Fix $p\in\left\{  1,2,\ldots,k\right\}  $. For each $i\in\mathbb{N}$, we let
$\mathcal{A}\left(  i\right)  $ be the statement%
\begin{equation}
\left(  \text{if }p+i\in\left\{  1,2,\ldots,k\right\}  \text{, then }%
a_{p}\equiv a_{p+i}\operatorname{mod}n\right)  . \label{pf.prop.mod.chain.Ai}%
\end{equation}

We shall prove that this statement $\mathcal{A}\left(  i\right)  $ holds for
each $i\in\mathbb{N}$.

In fact, let us prove this by induction on $i$:\ \ \ \ \footnote{Thus, the
letter \textquotedblleft$i$\textquotedblright\ plays the role of the
\textquotedblleft$n$\textquotedblright\ in Theorem \ref{thm.ind.IP0} (since we
are already using \textquotedblleft$n$\textquotedblright\ for a different
thing).}

\textit{Induction base:} The statement $\mathcal{A}\left(  0\right)  $
holds\footnote{\textit{Proof.} Proposition \ref{prop.mod.transi} \textbf{(a)}
(applied to $a=a_{p}$) yields $a_{p}\equiv a_{p}\operatorname{mod}n$. In view
of $p=p+0$, this rewrites as $a_{p}\equiv a_{p+0}\operatorname{mod}n$. Hence,
$\left(  \text{if }p+0\in\left\{  1,2,\ldots,k\right\}  \text{, then }%
a_{p}\equiv a_{p+0}\operatorname{mod}n\right)  $. But this is precisely the
statement $\mathcal{A}\left(  0\right)  $. Hence, the statement $\mathcal{A}%
\left(  0\right)  $ holds.}. This completes the induction base.

\textit{Induction step:} Let $m\in\mathbb{N}$. Assume that $\mathcal{A}\left(
m\right)  $ holds. We must show that $\mathcal{A}\left(  m+1\right)  $ holds.

We have assumed that $\mathcal{A}\left(  m\right)  $ holds. In other words,%
\begin{equation}
\left(  \text{if }p+m\in\left\{  1,2,\ldots,k\right\}  \text{, then }%
a_{p}\equiv a_{p+m}\operatorname{mod}n\right)  .
\label{pf.prop.mod.chain.Ai.IH}%
\end{equation}

Next, let us assume that $p+\left(  m+1\right)  \in\left\{  1,2,\ldots
,k\right\}  $. Thus, $p+\left(  m+1\right)  \leq k$, so that $p+m+1=p+\left(
m+1\right)  \leq k$ and therefore $p+m\leq k-1$. Also, $p\in\left\{
1,2,\ldots,k\right\}  $, so that $p\geq1$ and thus $\underbrace{p}_{\geq
1}+\underbrace{m}_{\geq0}\geq1+0=1$. Combining this with $p+m\leq k-1$, we
obtain $p+m\in\left\{  1,2,\ldots,k-1\right\}  \subseteq\left\{
1,2,\ldots,k\right\}  $. Hence, (\ref{pf.prop.mod.chain.Ai.IH}) shows that
$a_{p}\equiv a_{p+m}\operatorname{mod}n$. But (\ref{pf.prop.mod.chain.ass})
(applied to $p+m$ instead of $i$) yields $a_{p+m}\equiv a_{\left(  p+m\right)
+1}\operatorname{mod}n$ (since $p+m\in\left\{  1,2,\ldots,k-1\right\}  $).

So we know that $a_{p}\equiv a_{p+m}\operatorname{mod}n$ and $a_{p+m}\equiv
a_{\left(  p+m\right)  +1}\operatorname{mod}n$. Hence, Proposition
\ref{prop.mod.transi} \textbf{(c)} (applied to $a=a_{p}$, $b=a_{p+m}$ and
$c=a_{\left(  p+m\right)  +1}$) yields $a_{p}\equiv a_{\left(  p+m\right)
+1}\operatorname{mod}n$. Since $\left(  p+m\right)  +1=p+\left(  m+1\right)
$, this rewrites as $a_{p}\equiv a_{p+\left(  m+1\right)  }\operatorname{mod}%
n$.

Now, forget that we assumed that $p+\left(  m+1\right)  \in\left\{
1,2,\ldots,k\right\}  $. We thus have shown that%
\[
\left(  \text{if }p+\left(  m+1\right)  \in\left\{  1,2,\ldots,k\right\}
\text{, then }a_{p}\equiv a_{p+\left(  m+1\right)  }\operatorname{mod}%
n\right)  .
\]
But this is precisely the statement $\mathcal{A}\left(  m+1\right)  $. Thus,
$\mathcal{A}\left(  m+1\right)  $ holds.

Now, forget that we fixed $m$. We thus have shown that if $m\in\mathbb{N}$ is
such that $\mathcal{A}\left(  m\right)  $ holds, then $\mathcal{A}\left(
m+1\right)  $ also holds. This completes the induction step.

Thus, we have completed both the induction base and the induction step. Hence,
by induction, we conclude that $\mathcal{A}\left(  i\right)  $ holds for each
$i\in\mathbb{N}$. In other words, (\ref{pf.prop.mod.chain.Ai}) holds for each
$i\in\mathbb{N}$.

We are not done yet, since our goal is to prove Proposition
\ref{prop.mod.chain}, not merely to prove $\mathcal{A}\left(  i\right)  $. But
this is now easy.

First, let us forget that we fixed $p$. Thus, we have shown that
(\ref{pf.prop.mod.chain.Ai}) holds for each $p\in\left\{  1,2,\ldots
,k\right\}  $ and $i\in\mathbb{N}$.

But we have either $u\leq v$ or $u>v$. In other words, we are in one of the
following two cases:

\textit{Case 1:} We have $u\leq v$.

\textit{Case 2:} We have $u>v$.

Let us first consider Case 1. In this case, we have $u\leq v$. Thus,
$v-u\geq0$, so that $v-u\in\mathbb{N}$. But recall that
(\ref{pf.prop.mod.chain.Ai}) holds for each $p\in\left\{  1,2,\ldots
,k\right\}  $ and $i\in\mathbb{N}$. Applying this to $p=u$ and $i=v-u$, we
conclude that (\ref{pf.prop.mod.chain.Ai}) holds for $p=u$ and $i=v-u$ (since
$u\in\left\{  1,2,\ldots,k\right\}  $ and $v-u\in\mathbb{N}$). In other words,%
\[
\left(  \text{if }u+\left(  v-u\right)  \in\left\{  1,2,\ldots,k\right\}
\text{, then }a_{u}\equiv a_{u+\left(  v-u\right)  }\operatorname{mod}%
n\right)  .
\]
Since $u+\left(  v-u\right)  =v$, this rewrites as%
\[
\left(  \text{if }v\in\left\{  1,2,\ldots,k\right\}  \text{, then }a_{u}\equiv
a_{v}\operatorname{mod}n\right)  .
\]
Since $v\in\left\{  1,2,\ldots,k\right\}  $ holds (by assumption), we conclude
that $a_{u}\equiv a_{v}\operatorname{mod}n$. Thus, Proposition
\ref{prop.mod.chain} is proven in Case 1.

Let us now consider Case 2. In this case, we have $u>v$. Thus, $u-v>0$, so
that $u-v\in\mathbb{N}$. But recall that (\ref{pf.prop.mod.chain.Ai}) holds
for each $p\in\left\{  1,2,\ldots,k\right\}  $ and $i\in\mathbb{N}$. Applying
this to $p=v$ and $i=u-v$, we conclude that (\ref{pf.prop.mod.chain.Ai}) holds
for $p=v$ and $i=u-v$ (since $v\in\left\{  1,2,\ldots,k\right\}  $ and
$u-v\in\mathbb{N}$). In other words,%
\[
\left(  \text{if }v+\left(  u-v\right)  \in\left\{  1,2,\ldots,k\right\}
\text{, then }a_{v}\equiv a_{v+\left(  u-v\right)  }\operatorname{mod}%
n\right)  .
\]
Since $v+\left(  u-v\right)  =u$, this rewrites as%
\[
\left(  \text{if }u\in\left\{  1,2,\ldots,k\right\}  \text{, then }a_{v}\equiv
a_{u}\operatorname{mod}n\right)  .
\]
Since $u\in\left\{  1,2,\ldots,k\right\}  $ holds (by assumption), we conclude
that $a_{v}\equiv a_{u}\operatorname{mod}n$. Therefore, Proposition
\ref{prop.mod.transi} \textbf{(b)} (applied to $a=a_{v}$ and $b=a_{u}$) yields
that $a_{u}\equiv a_{v}\operatorname{mod}n$. Thus, Proposition
\ref{prop.mod.chain} is proven in Case 2.

Hence, Proposition \ref{prop.mod.chain} is proven in both Cases 1 and 2. Since
these two Cases cover all possibilities, we thus conclude that Proposition
\ref{prop.mod.chain} always holds.
\end{proof}

\subsubsection{Chains of inequalities (a digression)}

Much of the above proof of Proposition \ref{prop.mod.chain} was unremarkable
and straightforward reasoning -- but this proof is nevertheless fundamental
and important. More or less the same argument can be used to show the
following fact about chains of inequalities:

\begin{proposition}
\label{prop.mod.chain-ineq}Let $a_{1},a_{2},\ldots,a_{k}$ be integers such
that $a_{1}\leq a_{2}\leq\cdots\leq a_{k}$. (Recall that the statement
\textquotedblleft$a_{1}\leq a_{2}\leq\cdots\leq a_{k}$\textquotedblright%
\ means that $\left(  a_{i}\leq a_{i+1}\text{ holds for each }i\in\left\{
1,2,\ldots,k-1\right\}  \right)  $.) Let $u$ and $v$ be two elements of
$\left\{  1,2,\ldots,k\right\}  $ such that $u\leq v$. Then,%
\[
a_{u}\leq a_{v}.
\]

\end{proposition}

Proposition \ref{prop.mod.chain-ineq} is similar to Proposition
\ref{prop.mod.chain}, with the congruences replaced by inequalities; but note
that the condition \textquotedblleft$u\leq v$\textquotedblright\ is now
required. Make sure you understand where you need this condition when adapting
the proof of Proposition \ref{prop.mod.chain} to Proposition
\ref{prop.mod.chain-ineq}!

For future use, let us prove a corollary of Proposition
\ref{prop.mod.chain-ineq} which essentially observes that the inequality sign
in $a_{u}\leq a_{v}$ can be made strict if there is any strict inequality sign
between $a_{u}$ and $a_{v}$ in the chain $a_{1}\leq a_{2}\leq\cdots\leq a_{k}$:

\begin{corollary}
\label{cor.mod.chain-ineq2}Let $a_{1},a_{2},\ldots,a_{k}$ be integers such
that $a_{1}\leq a_{2}\leq\cdots\leq a_{k}$. Let $u$ and $v$ be two elements of
$\left\{  1,2,\ldots,k\right\}  $ such that $u\leq v$. Let $p\in\left\{
u,u+1,\ldots,v-1\right\}  $ be such that $a_{p}<a_{p+1}$. Then,%
\[
a_{u}<a_{v}.
\]

\end{corollary}

\begin{proof}
[Proof of Corollary \ref{cor.mod.chain-ineq2}.]From $u\in\left\{
1,2,\ldots,k\right\}  $, we obtain $u\geq1$. From $v\in\left\{  1,2,\ldots
,k\right\}  $, we obtain $v\leq k$. From $p\in\left\{  u,u+1,\ldots
,v-1\right\}  $, we obtain $p\geq u$ and $p\leq v-1$. From $p\leq v-1$, we
obtain $p+1\leq v\leq k$. Combining this with $p+1\geq p\geq u\geq1$, we
obtain $p+1\in\left\{  1,2,\ldots,k\right\}  $ (since $p+1$ is an integer).
Combining $p\leq v-1\leq v\leq k$ with $p\geq u\geq1$, we obtain $p\in\left\{
1,2,\ldots,k\right\}  $ (since $p$ is an integer). We thus know that both $p$
and $p+1$ are elements of $\left\{  1,2,\ldots,k\right\}  $.

We have $p\geq u$, thus $u\leq p$. Hence, Proposition
\ref{prop.mod.chain-ineq} (applied to $p$ instead of $v$) yields $a_{u}\leq
a_{p}$. Combining this with $a_{p}<a_{p+1}$, we find $a_{u}<a_{p+1}$.

We have $p+1\leq v$. Hence, Proposition \ref{prop.mod.chain-ineq} (applied to
$p+1$ instead of $u$) yields $a_{p+1}\leq a_{v}$. Combining $a_{u}<a_{p+1}$
with $a_{p+1}\leq a_{v}$, we obtain $a_{u}<a_{v}$. This proves Corollary
\ref{cor.mod.chain-ineq2}.
\end{proof}

In particular, we see that the inequality sign in $a_{u}\leq a_{v}$ is strict
when $u<v$ holds and \textbf{all} inequality signs in the chain $a_{1}\leq
a_{2}\leq\cdots\leq a_{k}$ are strict:

\begin{corollary}
\label{cor.mod.chain-ineq3}Let $a_{1},a_{2},\ldots,a_{k}$ be integers such
that $a_{1}<a_{2}<\cdots<a_{k}$. (Recall that the statement \textquotedblleft%
$a_{1}<a_{2}<\cdots<a_{k}$\textquotedblright\ means that $\left(
a_{i}<a_{i+1}\text{ holds for each }i\in\left\{  1,2,\ldots,k-1\right\}
\right)  $.) Let $u$ and $v$ be two elements of $\left\{  1,2,\ldots
,k\right\}  $ such that $u<v$. Then,%
\[
a_{u}<a_{v}.
\]

\end{corollary}

\begin{proof}
[Proof of Corollary \ref{cor.mod.chain-ineq3}.]From $u<v$, we obtain $u\leq
v-1$ (since $u$ and $v$ are integers). Combining this with $u\geq u$, we
conclude that $u\in\left\{  u,u+1,\ldots,v-1\right\}  $. Also, from
$a_{1}<a_{2}<\cdots<a_{k}$, we obtain $a_{1}\leq a_{2}\leq\cdots\leq a_{k}$.

We have $u\leq v-1$, thus $u+1\leq v\leq k$ (since $v\in\left\{
1,2,\ldots,k\right\}  $), so that $u\leq k-1$. Combining this with $u\geq1$
(which is a consequence of $u\in\left\{  1,2,\ldots,k\right\}  $), we find
$u\in\left\{  1,2,\ldots,k-1\right\}  $. Hence, from $a_{1}<a_{2}<\cdots
<a_{k}$, we obtain $a_{u}<a_{u+1}$. Hence, Corollary \ref{cor.mod.chain-ineq2}
(applied to $p=u$) yields $a_{u}<a_{v}$ (since $u\leq v$ (because $u<v$)).
This proves Corollary \ref{cor.mod.chain-ineq3}.
\end{proof}

\subsubsection{Addition, subtraction and multiplication of congruences}

Let us now return to the topic of congruences.

Chains of congruences can include equality signs. For example, if $a,b,c,d,n$
are integers, then \textquotedblleft$a\equiv b=c\equiv d\operatorname{mod}%
n$\textquotedblright\ means that $a\equiv b\operatorname{mod}n$ and $b=c$ and
$c\equiv d\operatorname{mod}n$. Such a chain is still a chain of congruences,
because $b=c$ implies $b\equiv c\operatorname{mod}n$ (by Proposition
\ref{prop.mod.transi} \textbf{(a)}).

Let us continue with basic properties of congruences:

\begin{proposition}
\label{prop.mod.+-*}Let $a$, $b$, $c$, $d$ and $n$ be integers such that
$a\equiv b\operatorname{mod}n$ and $c\equiv d\operatorname{mod}n$. Then:

\textbf{(a)} We have $a+c\equiv b+d\operatorname{mod}n$.

\textbf{(b)} We have $a-c\equiv b-d\operatorname{mod}n$.

\textbf{(c)} We have $ac\equiv bd\operatorname{mod}n$.
\end{proposition}

Note that Proposition \ref{prop.mod.+-*} does \textbf{not} claim that
$a/c\equiv b/d\operatorname{mod}n$. Indeed, this would not be true in general.
One reason for this is that $a/c$ and $b/d$ aren't always integers. But even
when they are, they may not satisfy $a/c\equiv b/d\operatorname{mod}n$. For
example, $6\equiv4\operatorname{mod}2$ and $2\equiv2\operatorname{mod}2$, but
$6/2\not \equiv 4/2\operatorname{mod}2$. Likewise, Proposition
\ref{prop.mod.+-*} does \textbf{not} claim that $a^{c}\equiv b^{d}%
\operatorname{mod}n$ even when $a,b,c,d$ are nonnegative; that too would not
be true. But we will soon see that a weaker statement (Proposition
\ref{prop.mod.pow}) holds. First, let us prove Proposition \ref{prop.mod.+-*}:

\begin{proof}
[Proof of Proposition \ref{prop.mod.+-*}.]From $a\equiv b\operatorname{mod}n$,
we conclude that $a$ is congruent to $b$ modulo $n$. In other words, $n\mid
a-b$ (by the definition of \textquotedblleft congruent\textquotedblright). In
other words, $n$ divides $a-b$. In other words, there exists an integer $w$
such that $a-b=nw$ (by the definition of \textquotedblleft
divides\textquotedblright). Consider this $w$, and denote it by $q$. Thus, $q$
is an integer such that $a-b=nq$.

\begin{vershort}
Similarly, from $c\equiv d\operatorname{mod}n$, we can construct an integer
$r$ such that $c-d=nr$. Consider this $r$.
\end{vershort}

\begin{verlong}
From $c\equiv d\operatorname{mod}n$, we conclude that $c$ is congruent to $d$
modulo $n$. In other words, $n\mid c-d$ (by the definition of
\textquotedblleft congruent\textquotedblright). In other words, $n$ divides
$c-d$. In other words, there exists an integer $w$ such that $c-d=nw$ (by the
definition of \textquotedblleft divides\textquotedblright). Consider this $w$,
and denote it by $r$. Thus, $r$ is an integer such that $c-d=nr$.
\end{verlong}

\textbf{(a)} We have%
\[
\left(  a+c\right)  -\left(  b+d\right)  =\underbrace{\left(  a-b\right)
}_{=nq}+\underbrace{\left(  c-d\right)  }_{=nr}=nq+nr=n\left(  q+r\right)  .
\]
Hence, there exists an integer $w$ such that $\left(  a+c\right)  -\left(
b+d\right)  =nw$ (namely, $w=q+r$). In other words, $n$ divides $\left(
a+c\right)  -\left(  b+d\right)  $ (by the definition of \textquotedblleft
divides\textquotedblright). In other words, $n\mid\left(  a+c\right)  -\left(
b+d\right)  $. In other words, $a+c\equiv b+d\operatorname{mod}n$ (by the
definition of \textquotedblleft congruent\textquotedblright). This proves
Proposition \ref{prop.mod.+-*} \textbf{(a)}.

\textbf{(b)} We have%
\[
\left(  a-c\right)  -\left(  b-d\right)  =\underbrace{\left(  a-b\right)
}_{=nq}-\underbrace{\left(  c-d\right)  }_{=nr}=nq-nr=n\left(  q-r\right)  .
\]
Hence, there exists an integer $w$ such that $\left(  a-c\right)  -\left(
b-d\right)  =nw$ (namely, $w=q-r$). In other words, $n$ divides $\left(
a-c\right)  -\left(  b-d\right)  $ (by the definition of \textquotedblleft
divides\textquotedblright). In other words, $n\mid\left(  a-c\right)  -\left(
b-d\right)  $. In other words, $a-c\equiv b-d\operatorname{mod}n$ (by the
definition of \textquotedblleft congruent\textquotedblright). This proves
Proposition \ref{prop.mod.+-*} \textbf{(b)}.

\textbf{(c)} We have $ac-ad=a\underbrace{\left(  c-d\right)  }_{=nr}%
=anr=n\left(  ar\right)  $. Hence, there exists an integer $w$ such that
$ac-ad=nw$ (namely, $w=ar$). In other words, $n$ divides $ac-ad$ (by the
definition of \textquotedblleft divides\textquotedblright). In other words,
$n\mid ac-ad$. In other words, $ac\equiv ad\operatorname{mod}n$ (by the
definition of \textquotedblleft congruent\textquotedblright).

We have $ad-bd=\underbrace{\left(  a-b\right)  }_{=nq}d=nqd=n\left(
qd\right)  $. Hence, there exists an integer $w$ such that $ad-bd=nw$ (namely,
$w=qd$). In other words, $n$ divides $ad-bd$ (by the definition of
\textquotedblleft divides\textquotedblright). In other words, $n\mid ad-bd$.
In other words, $ad\equiv bd\operatorname{mod}n$ (by the definition of
\textquotedblleft congruent\textquotedblright).

Now, we know that $ac\equiv ad\operatorname{mod}n$ and $ad\equiv
bd\operatorname{mod}n$. Hence, Proposition \ref{prop.mod.transi} \textbf{(c)}
(applied to $ac$, $ad$ and $bd$ instead of $a$, $b$ and $c$) shows that
$ac\equiv bd\operatorname{mod}n$. This proves Proposition \ref{prop.mod.+-*}
\textbf{(c)}.
\end{proof}

Proposition \ref{prop.mod.+-*} shows yet another aspect in which congruences
(modulo $n$) behave like equalities: They can be added, subtracted and
multiplied, in the following sense:

\begin{itemize}
\item We can add two congruences modulo $n$ (in the sense of adding each side
of one congruence to the corresponding side of the other); this yields a new
congruence modulo $n$ (because of Proposition \ref{prop.mod.+-*} \textbf{(a)}).

\item We can subtract two congruences modulo $n$; this yields a new congruence
modulo $n$ (because of Proposition \ref{prop.mod.+-*} \textbf{(b)}).

\item We can multiply two congruences modulo $n$; this yields a new congruence
modulo $n$ (because of Proposition \ref{prop.mod.+-*} \textbf{(c)}).
\end{itemize}

\subsubsection{Substitutivity for congruences}

Combined with Proposition \ref{prop.mod.transi}, these observations lead to a
further feature of congruences, which is even more important: the principle of
\textit{substitutivity for congruences}. We are not going to state it fully
formally (as it is a meta-mathematical principle), but merely explain what it means.

Recall that the \textit{principle of substitutivity for equalities} says the following:

\begin{statement}
\textit{Principle of substitutivity for equalities:} If two
objects\footnote{\textquotedblleft Objects\textquotedblright\ can be numbers,
sets, tuples or any other mathematical objects.} $x$ and $x^{\prime}$ are
equal, and if we have any expression $A$ that involves the object $x$, then we
can replace this $x$ (or, more precisely, any arbitrary appearance of $x$ in
$A$) in $A$ by $x^{\prime}$; the value of the resulting expression $A^{\prime
}$ will be equal to the value of $A$.
\end{statement}

Here are two examples of how this principle can be used:

\begin{itemize}
\item If $a,b,c,d,e,c^{\prime}$ are numbers such that $c=c^{\prime}$, then the
principle of substitutivity for equalities says that we can replace $c$ by
$c^{\prime}$ in the expression $a\left(  b-\left(  c+d\right)  e\right)  $,
and the value of the resulting expression $a\left(  b-\left(  c^{\prime
}+d\right)  e\right)  $ will be equal to the value of $a\left(  b-\left(
c+d\right)  e\right)  $; that is, we have%
\begin{equation}
a\left(  b-\left(  c+d\right)  e\right)  =a\left(  b-\left(  c^{\prime
}+d\right)  e\right)  . \label{eq.mod.substitutivity-nums.1}%
\end{equation}

\item If $a,b,c,a^{\prime}$ are numbers such that $a=a^{\prime}$, then
\begin{equation}
\left(  a-b\right)  \left(  a+b\right)  =\left(  a^{\prime}-b\right)  \left(
a+b\right)  , \label{eq.mod.substitutivity-nums.2}%
\end{equation}
because the principle of substitutivity allows us to replace the first $a$
appearing in the expression $\left(  a-b\right)  \left(  a+b\right)  $ by an
$a^{\prime}$. (We can also replace the second $a$ by $a^{\prime}$, of course.)
\end{itemize}

More generally, we can make several such replacements at the same time.

The principle of substitutivity for equalities is one of the headstones of
mathematical logic; it is the essence of what it means for two objects to be equal.

The \textit{principle of substitutivity for congruences} is similar, but far
less fundamental; it says the following:

\begin{statement}
\textit{Principle of substitutivity for congruences:} Fix an integer $n$. If
two numbers $x$ and $x^{\prime}$ are congruent to each other modulo $n$ (that
is, $x\equiv x^{\prime}\operatorname{mod}n$), and if we have any expression
$A$ that involves only integers, addition, subtraction and multiplication, and
involves the object $x$, then we can replace this $x$ (or, more precisely, any
arbitrary appearance of $x$ in $A$) in $A$ by $x^{\prime}$; the value of the
resulting expression $A^{\prime}$ will be congruent to the value of $A$ modulo
$n$.
\end{statement}

Note that this principle is less general than the principle of substitutivity
for equalities, because it only applies to expressions that are built from
integers and certain operations (note that division is not one of these
operations). But it still lets us prove analogues of our above examples
(\ref{eq.mod.substitutivity-nums.1}) and (\ref{eq.mod.substitutivity-nums.2}):

\begin{itemize}
\item If $n$ is any integer, and if $a,b,c,d,e,c^{\prime}$ are integers such
that $c\equiv c^{\prime}\operatorname{mod}n$, then the principle of
substitutivity for congruences says that we can replace $c$ by $c^{\prime}$ in
the expression $a\left(  b-\left(  c+d\right)  e\right)  $, and the value of
the resulting expression $a\left(  b-\left(  c^{\prime}+d\right)  e\right)  $
will be congruent to the value of $a\left(  b-\left(  c+d\right)  e\right)  $
modulo $n$; that is, we have%
\begin{equation}
a\left(  b-\left(  c+d\right)  e\right)  \equiv a\left(  b-\left(  c^{\prime
}+d\right)  e\right)  \operatorname{mod}n.
\label{eq.mod.substitutivity-congs.1}%
\end{equation}

\item If $n$ is any integer, and if $a,b,c,a^{\prime}$ are integers such that
$a\equiv a^{\prime}\operatorname{mod}n$, then
\begin{equation}
\left(  a-b\right)  \left(  a+b\right)  \equiv\left(  a^{\prime}-b\right)
\left(  a+b\right)  \operatorname{mod}n, \label{eq.mod.substitutivity-congs.2}%
\end{equation}
because the principle of substitutivity allows us to replace the first $a$
appearing in the expression $\left(  a-b\right)  \left(  a+b\right)  $ by an
$a^{\prime}$. (We can also replace the second $a$ by $a^{\prime}$, of course.)
\end{itemize}

We shall not prove the principle of substitutivity for congruences, since we
have not formalized it (after all, we have not defined what an
\textquotedblleft expression\textquotedblright\ is). But we shall prove the
specific congruences (\ref{eq.mod.substitutivity-congs.1}) and
(\ref{eq.mod.substitutivity-congs.2}) using Proposition \ref{prop.mod.+-*} and
Proposition \ref{prop.mod.transi}; the way in which we prove these congruences
is symptomatic: Every congruence obtained from the principle of substitutivity
for congruences can be proven in a manner like these. Thus, we hope that the
proofs of (\ref{eq.mod.substitutivity-congs.1}) and
(\ref{eq.mod.substitutivity-congs.2}) given below serve as templates which can
easily be adapted to any other situation in which an application of the
principle of substitutivity for congruences needs to be justified.

\begin{proof}
[Proof of (\ref{eq.mod.substitutivity-congs.1}).]Let $n$ be any integer, and
let $a,b,c,d,e,c^{\prime}$ be integers such that $c\equiv c^{\prime
}\operatorname{mod}n$.

Adding the congruence $c\equiv c^{\prime}\operatorname{mod}n$ with the
congruence $d\equiv d\operatorname{mod}n$ (which follows from Proposition
\ref{prop.mod.transi} \textbf{(a)}), we obtain $c+d\equiv c^{\prime
}+d\operatorname{mod}n$. Multiplying this congruence with the congruence
$e\equiv e\operatorname{mod}n$ (which follows from Proposition
\ref{prop.mod.transi} \textbf{(a)}), we obtain $\left(  c+d\right)
e\equiv\left(  c^{\prime}+d\right)  e\operatorname{mod}n$. Subtracting this
congruence from the congruence $b\equiv b\operatorname{mod}n$ (which, again,
follows from Proposition \ref{prop.mod.transi} \textbf{(a)}), we obtain
$b-\left(  c+d\right)  e\equiv b-\left(  c^{\prime}+d\right)
e\operatorname{mod}n$. Multiplying the congruence $a\equiv a\operatorname{mod}%
n$ (which follows from Proposition \ref{prop.mod.transi} \textbf{(a)}) with
this congruence, we obtain $a\left(  b-\left(  c+d\right)  e\right)  \equiv
a\left(  b-\left(  c^{\prime}+d\right)  e\right)  \operatorname{mod}n$. This
proves (\ref{eq.mod.substitutivity-congs.1}).
\end{proof}

\begin{proof}
[Proof of (\ref{eq.mod.substitutivity-congs.2}).]Let $n$ be any integer, and
let $a,b,c,a^{\prime}$ be integers such that $a\equiv a^{\prime}%
\operatorname{mod}n$.

Subtracting the congruence $b\equiv b\operatorname{mod}n$ (which follows from
Proposition \ref{prop.mod.transi} \textbf{(a)}) from the congruence $a\equiv
a^{\prime}\operatorname{mod}n$, we obtain $a-b\equiv a^{\prime}%
-b\operatorname{mod}n$. Multiplying this congruence with the congruence
$a+b\equiv a+b\operatorname{mod}n$ (which follows from Proposition
\ref{prop.mod.transi} \textbf{(a)}), we obtain $\left(  a-b\right)  \left(
a+b\right)  \equiv\left(  a^{\prime}-b\right)  \left(  a+b\right)
\operatorname{mod}n$. This proves (\ref{eq.mod.substitutivity-congs.2}).
\end{proof}

As we said, these two proofs are exemplary: Any congruence obtained from the
principle of substitutivity for congruences can be proven in such a way
(starting with the congruence $x\equiv x^{\prime}\operatorname{mod}n$, and
then \textquotedblleft wrapping\textquotedblright\ it up in the expression $A$
by repeatedly adding, multiplying and subtracting congruences that follow from
Proposition \ref{prop.mod.transi} \textbf{(a)}).

When we apply the principle of substitutivity for congruences, we shall use
underbraces to point out which integers we are replacing. For example, when
deriving (\ref{eq.mod.substitutivity-congs.1}) from this principle, we shall
write%
\[
a\left(  b-\left(  \underbrace{c}_{\equiv c^{\prime}\operatorname{mod}%
n}+d\right)  e\right)  \equiv a\left(  b-\left(  c^{\prime}+d\right)
e\right)  \operatorname{mod}n,
\]
in order to stress that we are replacing $c$ by $c^{\prime}$. Likewise, when
deriving (\ref{eq.mod.substitutivity-congs.2}) from this principle, we shall
write%
\[
\left(  \underbrace{a}_{\equiv a^{\prime}\operatorname{mod}n}-b\right)
\left(  a+b\right)  \equiv\left(  a^{\prime}-b\right)  \left(  a+b\right)
\operatorname{mod}n,
\]
in order to stress that we are replacing the first $a$ (but not the second
$a$) by $a^{\prime}$.

The principle of substitutivity for congruences allows us to replace a
\textbf{single} integer $x$ appearing in an expression by another integer
$x^{\prime}$ that is congruent to $x$ modulo $n$. Applying this principle many
times, we thus conclude that we can also replace \textbf{several} integers at
the same time (because we can get to the same result by performing these
replacements one at a time, and Proposition \ref{prop.mod.chain} shows that
the final result will be congruent to the original result).

For example, if seven integers $a,a^{\prime},b,b^{\prime},c,c^{\prime},n$
satisfy $a\equiv a^{\prime}\operatorname{mod}n$ and $b\equiv b^{\prime
}\operatorname{mod}n$ and $c\equiv c^{\prime}\operatorname{mod}n$, then%
\begin{equation}
bc+ca+ab\equiv b^{\prime}c^{\prime}+c^{\prime}a^{\prime}+a^{\prime}b^{\prime
}\operatorname{mod}n, \label{eq.mod.substitutivity-congs.3}%
\end{equation}
because we can replace all the six integers $b,c,c,a,a,b$ in the expression
$bc+ca+ab$ (listed in the order of their appearance in this expression) by
$b^{\prime},c^{\prime},c^{\prime},a^{\prime},a^{\prime},b^{\prime}$,
respectively. If we want to derive this from the principle of substitutivity
for congruences, we must perform the replacements one at a time, e.g., as
follows:%
\begin{align*}
\underbrace{b}_{\equiv b^{\prime}\operatorname{mod}n}c+ca+ab  &  \equiv
b^{\prime}\underbrace{c}_{\equiv c^{\prime}\operatorname{mod}n}+ca+ab\equiv
b^{\prime}c^{\prime}+\underbrace{c}_{\equiv c^{\prime}\operatorname{mod}%
n}a+ab\\
&  \equiv b^{\prime}c^{\prime}+c^{\prime}\underbrace{a}_{\equiv a^{\prime
}\operatorname{mod}n}+ab\equiv b^{\prime}c^{\prime}+c^{\prime}a^{\prime
}+\underbrace{a}_{\equiv a^{\prime}\operatorname{mod}n}b\\
&  \equiv b^{\prime}c^{\prime}+c^{\prime}a^{\prime}+a^{\prime}\underbrace{b}%
_{\equiv b^{\prime}\operatorname{mod}n}\equiv b^{\prime}c^{\prime}+c^{\prime
}a^{\prime}+a^{\prime}b^{\prime}\operatorname{mod}n.
\end{align*}
Of course, we shall always just show the replacements as a single step:%
\[
\underbrace{b}_{\equiv b^{\prime}\operatorname{mod}n}\underbrace{c}_{\equiv
c^{\prime}\operatorname{mod}n}+\underbrace{c}_{\equiv c^{\prime}%
\operatorname{mod}n}\underbrace{a}_{\equiv a^{\prime}\operatorname{mod}%
n}+\underbrace{a}_{\equiv a^{\prime}\operatorname{mod}n}\underbrace{b}_{\equiv
b^{\prime}\operatorname{mod}n}\equiv b^{\prime}c^{\prime}+c^{\prime}a^{\prime
}+a^{\prime}b^{\prime}\operatorname{mod}n.
\]

\begin{noncompile}
(because multiplying the congruences $b\equiv b^{\prime}\operatorname{mod}n$
and $c\equiv c^{\prime}\operatorname{mod}n$ yields $bc\equiv b^{\prime
}c^{\prime}\operatorname{mod}n$, and similarly we can show that $ca\equiv
c^{\prime}a^{\prime}\operatorname{mod}n$ and $ab\equiv a^{\prime}b^{\prime
}\operatorname{mod}n$, and then we can add the last three congruences together
to obtain (\ref{eq.mod.add-mult.ex1})). This principle is called
\textit{substitutivity for congruences}. When we apply this principle, we
shall use underbraces to point out which congruences we are adding,
subtracting or multiplying. For example, the proof of
(\ref{eq.mod.substitutivity-congs.3}) is thus written as%
\[
\underbrace{b}_{\equiv b^{\prime}\operatorname{mod}n}\underbrace{c}_{\equiv
c^{\prime}\operatorname{mod}n}+\underbrace{c}_{\equiv c^{\prime}%
\operatorname{mod}n}\underbrace{a}_{\equiv a^{\prime}\operatorname{mod}%
n}+\underbrace{a}_{\equiv a^{\prime}\operatorname{mod}n}\underbrace{b}_{\equiv
b^{\prime}\operatorname{mod}n}\equiv b^{\prime}c^{\prime}+c^{\prime}a^{\prime
}+a^{\prime}b^{\prime}\operatorname{mod}n.
\]
Note that Proposition \ref{prop.mod.transi} \textbf{(a)} shows that an integer
is always congruent to itself modulo $n$ (whatever $n$ is); therefore, we can
replace an integer in a congruence by itself. Thus, when we apply
substitutivity for congruences, we don't need to replace \textbf{all}
integers. For example, if four integers $a,b,a^{\prime},n$ satisfy $a\equiv
a^{\prime}\operatorname{mod}n$, then%
\[
\left(  \underbrace{a}_{\equiv a^{\prime}\operatorname{mod}n}+b\right)
\left(  a^{\prime}+b\right)  \equiv\left(  a^{\prime}+b\right)  \left(
a^{\prime}+b\right)  =\left(  a^{\prime}+b\right)  ^{2}\operatorname{mod}n.
\]
\footnote{A way to prove this without referring to substitutivity is as
follows: Let $a,b,a^{\prime},n$ be four integers satisfying $a\equiv
a^{\prime}\operatorname{mod}n$. Then, adding the congruences $a\equiv
a^{\prime}\operatorname{mod}n$ and $b\equiv b\operatorname{mod}n$ (the second
of which follows from Proposition \ref{prop.mod.transi} \textbf{(a)}), we
obtain $a+b\equiv a^{\prime}+b\operatorname{mod}n$. Multiplying this
congruence with the congruence $a^{\prime}+b\equiv a^{\prime}%
+b\operatorname{mod}n$ (which follows from Proposition \ref{prop.mod.transi}
\textbf{(a)}), we obtain $\left(  a+b\right)  \left(  a^{\prime}+b\right)
\equiv\left(  a^{\prime}+b\right)  \left(  a^{\prime}+b\right)
\operatorname{mod}n$, which is what we wanted to prove. Any time you apply
substitutivity for congruences, you are implicitly making an argument like
this.}
\end{noncompile}

\subsubsection{Taking congruences to the $k$-th power}

We have seen that congruences (like equalities) can be added, subtracted and
multiplied (but, unlike equalities, they cannot be divided). One other thing
we can do with congruences is taking powers of them, as long as the exponent
is a nonnegative integer. This relies on the following fact:

\begin{proposition}
\label{prop.mod.pow}Let $a$, $b$ and $n$ be three integers such that $a\equiv
b\operatorname{mod}n$. Then, $a^{k}\equiv b^{k}\operatorname{mod}n$ for each
$k\in\mathbb{N}$.
\end{proposition}

The following proof of Proposition \ref{prop.mod.pow} is an example of a
straightforward inductive proof; the only thing to keep in mind is that it
uses induction on $k$, not induction on $n$ as some of our previous proofs did.

\begin{proof}
[Proof of Proposition \ref{prop.mod.pow}.]We claim that%
\begin{equation}
a^{k}\equiv b^{k}\operatorname{mod}n\ \ \ \ \ \ \ \ \ \ \text{for each }%
k\in\mathbb{N}\text{.} \label{pf.prop.mod.pow.goal}%
\end{equation}

We shall prove (\ref{pf.prop.mod.pow.goal}) by induction on $k$:

\textit{Induction base:} We have $1\equiv1\operatorname{mod}n$ (by Proposition
\ref{prop.mod.transi} \textbf{(a)}). In view of $a^{0}=1$ and $b^{0}=1$, this
rewrites as $a^{0}\equiv b^{0}\operatorname{mod}n$. In other words,
(\ref{pf.prop.mod.pow.goal}) holds for $k=0$. This completes the induction base.

\textit{Induction step:} Let $m\in\mathbb{N}$. Assume that
(\ref{pf.prop.mod.pow.goal}) holds for $k=m$. We must show that
(\ref{pf.prop.mod.pow.goal}) holds for $k=m+1$.

We have assumed that (\ref{pf.prop.mod.pow.goal}) holds for $k=m$. In other
words, we have $a^{m}\equiv b^{m}\operatorname{mod}n$. Now,%
\[
a^{m+1}=\underbrace{a^{m}}_{\equiv b^{m}\operatorname{mod}n}\underbrace{a}%
_{\equiv b\operatorname{mod}n}\equiv b^{m}b=b^{m+1}\operatorname{mod}n.
\]
\footnote{This computation relied on the principle of substitutivity for
congruences. Here is how to rewrite this argument in a more explicit way
(without using this principle): We have $a^{m}\equiv b^{m}\operatorname{mod}n$
and $a\equiv b\operatorname{mod}n$. Hence, Proposition \ref{prop.mod.+-*}
\textbf{(c)} (applied to $a^{m}$, $b^{m}$, $a$ and $b$ instead of $a$, $b$,
$c$ and $d$) yields $a^{m}a\equiv b^{m}b\operatorname{mod}n$. This rewrites as
$a^{m+1}\equiv b^{m+1}\operatorname{mod}n$ (since $a^{m+1}=a^{m}a$ and
$b^{m+1}=b^{m}b$).} In other words, (\ref{pf.prop.mod.pow.goal}) holds for
$k=m+1$. This completes the induction step. Hence, (\ref{pf.prop.mod.pow.goal}%
) is proven by induction. This proves Proposition \ref{prop.mod.pow}.
\end{proof}

\subsection{A few recursively defined sequences}

\subsubsection{$a_{n}=a_{n-1}^{q}+r$}

We next proceed to give some more examples of proofs by induction.

\begin{example}
\label{exa.rec-seq.1}Let $\left(  a_{0},a_{1},a_{2},\ldots\right)  $ be a
sequence of integers defined recursively by%
\begin{align*}
a_{0}  &  =0,\ \ \ \ \ \ \ \ \ \ \text{and}\\
a_{n}  &  =a_{n-1}^{2}+1\ \ \ \ \ \ \ \ \ \ \text{for each }n\geq1.
\end{align*}
(\textquotedblleft Defined recursively\textquotedblright\ means that we aren't
defining each entry $a_{n}$ of our sequence by an explicit formula, but rather
defining $a_{n}$ in terms of the previous entries $a_{0},a_{1},\ldots,a_{n-1}%
$. Thus, in order to compute some entry $a_{n}$ of our sequence, we need to
compute all the previous entries $a_{0},a_{1},\ldots,a_{n-1}$. This means that
if we want to compute $a_{n}$, we should first compute $a_{0}$, then compute
$a_{1}$ (using our value of $a_{0}$), then compute $a_{2}$ (using our values
of $a_{0}$ and $a_{1}$), and so on, until we reach $a_{n}$. For example, in
order to compute $a_{6}$, we proceed as follows:%
\begin{align*}
a_{0}  &  =0;\\
a_{1}  &  =a_{0}^{2}+1=0^{2}+1=1;\\
a_{2}  &  =a_{1}^{2}+1=1^{2}+1=2;\\
a_{3}  &  =a_{2}^{2}+1=2^{2}+1=5;\\
a_{4}  &  =a_{3}^{2}+1=5^{2}+1=26;\\
a_{5}  &  =a_{4}^{2}+1=26^{2}+1=677;\\
a_{6}  &  =a_{5}^{2}+1=677^{2}+1=458\,330.
\end{align*}
And similarly we can compute $a_{n}$ for any $n\in\mathbb{N}$.)

This sequence $\left(  a_{0},a_{1},a_{2},\ldots\right)  $ is not unknown: It
is the sequence \href{https://oeis.org/A003095}{A003095} in
\href{https://oeis.org/}{the Online Encyclopedia of Integer Sequences}.

A look at the first few entries of the sequence makes us realize that both
$a_{2}$ and $a_{3}$ divide $a_{6}$, just as the integers $2$ and $3$
themselves divide $6$. This suggests that we might have $a_{u}\mid a_{v}$
whenever $u$ and $v$ are two nonnegative integers satisfying $u\mid v$. We
shall soon prove this observation (which was found by Michael Somos in 2013)
in greater generality.
\end{example}

\begin{theorem}
\label{thm.rec-seq.somos-simple}Fix some $q\in\mathbb{N}$ and $r\in\mathbb{Z}%
$. Let $\left(  a_{0},a_{1},a_{2},\ldots\right)  $ be a sequence of integers
defined recursively by%
\begin{align*}
a_{0}  &  =0,\ \ \ \ \ \ \ \ \ \ \text{and}\\
a_{n}  &  =a_{n-1}^{q}+r\ \ \ \ \ \ \ \ \ \ \text{for each }n\geq1.
\end{align*}
(Note that if $q=2$ and $r=1$, then this sequence $\left(  a_{0},a_{1}%
,a_{2},\ldots\right)  $ is precisely the sequence $\left(  a_{0},a_{1}%
,a_{2},\ldots\right)  $ from Example \ref{exa.rec-seq.1}. If $q=3$ and $r=1$,
then our sequence $\left(  a_{0},a_{1},a_{2},\ldots\right)  $ is the sequence
\href{https://oeis.org/A135361}{A135361} in \href{https://oeis.org/}{the
Online Encyclopedia of Integer Sequences}. If $q=0$, then our sequence
$\left(  a_{0},a_{1},a_{2},\ldots\right)  $ is $\left(  0,r+1,r+1,r+1,\ldots
\right)  $. If $q=1$, then our sequence $\left(  a_{0},a_{1},a_{2}%
,\ldots\right)  $ is $\left(  0,r,2r,3r,4r,\ldots\right)  $, as can be easily
proven by induction.)

\textbf{(a)} For any $k\in\mathbb{N}$ and $n\in\mathbb{N}$, we have
$a_{k+n}\equiv a_{k}\operatorname{mod}a_{n}$.

\textbf{(b)} For any $n\in\mathbb{N}$ and $w\in\mathbb{N}$, we have $a_{n}\mid
a_{nw}$.

\textbf{(c)} If $u$ and $v$ are two nonnegative integers satisfying $u\mid v$,
then $a_{u}\mid a_{v}$.
\end{theorem}

\begin{proof}
[Proof of Theorem \ref{thm.rec-seq.somos-simple}.]\textbf{(a)} Let
$n\in\mathbb{N}$. We claim that%
\begin{equation}
a_{k+n}\equiv a_{k}\operatorname{mod}a_{n}\ \ \ \ \ \ \ \ \ \ \text{for every
}k\in\mathbb{N}. \label{pf.thm.rec-seq.somos-simple.a.claim}%
\end{equation}

We shall prove (\ref{pf.thm.rec-seq.somos-simple.a.claim}) by induction on $k$:

\textit{Induction base:} Proposition \ref{prop.mod.n=0} (applied to $a_{n}$
instead of $n$) yields $a_{n}\equiv0\operatorname{mod}a_{n}$. This rewrites as
$a_{n}\equiv a_{0}\operatorname{mod}a_{n}$ (since $a_{0}=0$). In other words,
$a_{0+n}\equiv a_{0}\operatorname{mod}a_{n}$ (since $0+n=n$). In other words,
(\ref{pf.thm.rec-seq.somos-simple.a.claim}) holds for $k=0$. This completes
the induction base.

\textit{Induction step:} Let $m\in\mathbb{N}$. Assume that
(\ref{pf.thm.rec-seq.somos-simple.a.claim}) holds for $k=m$. We must prove
that (\ref{pf.thm.rec-seq.somos-simple.a.claim}) holds for $k=m+1$.

We have assumed that (\ref{pf.thm.rec-seq.somos-simple.a.claim}) holds for
$k=m$. In other words, we have%
\[
a_{m+n}\equiv a_{m}\operatorname{mod}a_{n}.
\]
Hence, Proposition \ref{prop.mod.pow} (applied to $a_{m+n}$, $a_{m}$, $a_{n}$
and $q$ instead of $a$, $b$, $n$ and $k$) shows that $a_{m+n}^{q}\equiv
a_{m}^{q}\operatorname{mod}a_{n}$. Hence, $a_{m+n}^{q}+r\equiv a_{m}%
^{q}+r\operatorname{mod}a_{n}$. (Indeed, this follows by adding the congruence
$a_{m+n}^{q}\equiv a_{m}^{q}\operatorname{mod}a_{n}$ to the congruence
$r\equiv r\operatorname{mod}a_{n}$; the latter congruence is a consequence of
Proposition \ref{prop.mod.transi} \textbf{(a)}.)

Now, $\left(  m+1\right)  +n=\left(  m+n\right)  +1\geq1$. Hence, the
recursive definition of the sequence $\left(  a_{0},a_{1},a_{2},\ldots\right)
$ yields%
\[
a_{\left(  m+1\right)  +n}=a_{\left(  \left(  m+1\right)  +n\right)  -1}%
^{q}+r=a_{m+n}^{q}+r
\]
(since $\left(  \left(  m+1\right)  +n\right)  -1=m+n$). Also, $m+1\geq1$.
Hence, the recursive definition of the sequence $\left(  a_{0},a_{1}%
,a_{2},\ldots\right)  $ yields%
\[
a_{m+1}=a_{\left(  m+1\right)  -1}^{q}+r=a_{m}^{q}+r.
\]
The congruence $a_{m+n}^{q}+r\equiv a_{m}^{q}+r\operatorname{mod}a_{n}$
rewrites as $a_{\left(  m+1\right)  +n}\equiv a_{m+1}\operatorname{mod}a_{n}$
(since $a_{\left(  m+1\right)  +n}=a_{m+n}^{q}+r$ and $a_{m+1}=a_{m}^{q}+r$).
In other words, (\ref{pf.thm.rec-seq.somos-simple.a.claim}) holds for $k=m+1$.
This completes the induction step. Thus,
(\ref{pf.thm.rec-seq.somos-simple.a.claim}) is proven by induction.

Therefore, Theorem \ref{thm.rec-seq.somos-simple} \textbf{(a)} is proven.

\textbf{(b)} Let $n\in\mathbb{N}$. We claim that%
\begin{equation}
a_{n}\mid a_{nw}\ \ \ \ \ \ \ \ \ \ \text{for every }w\in\mathbb{N}.
\label{pf.thm.rec-seq.somos-simple.b.claim}%
\end{equation}

We shall prove (\ref{pf.thm.rec-seq.somos-simple.b.claim}) by induction on $w$:

\textit{Induction base:} We have $a_{n\cdot0}=a_{0}=0$. But $a_{n}\mid0$
(since $0=0a_{n}$). This rewrites as $a_{n}\mid a_{n\cdot0}$ (since
$a_{n\cdot0}=0$). In other words, (\ref{pf.thm.rec-seq.somos-simple.b.claim})
holds for $w=0$. This completes the induction base.

\textit{Induction step:} Let $m\in\mathbb{N}$. Assume that
(\ref{pf.thm.rec-seq.somos-simple.b.claim}) holds for $w=m$. We must prove
that (\ref{pf.thm.rec-seq.somos-simple.b.claim}) holds for $w=m+1$.

We have assumed that (\ref{pf.thm.rec-seq.somos-simple.b.claim}) holds for
$w=m$. In other words, we have $a_{n}\mid a_{nm}$.

Proposition \ref{prop.mod.0} \textbf{(a)} (applied to $a_{nm}$ and $a_{n}$
instead of $a$ and $n$) shows that we have $a_{nm}\equiv0\operatorname{mod}%
a_{n}$ if and only if $a_{n}\mid a_{nm}$. Hence, we have $a_{nm}%
\equiv0\operatorname{mod}a_{n}$ (since $a_{n}\mid a_{nm}$).

Theorem \ref{thm.rec-seq.somos-simple} \textbf{(a)} (applied to $k=nm$) yields
$a_{nm+n}\equiv a_{nm}\operatorname{mod}a_{n}$. Thus, $a_{nm+n}\equiv
a_{nm}\equiv0\operatorname{mod}a_{n}$. This is a chain of congruences; hence,
an application of Proposition \ref{prop.mod.chain} shows that $a_{nm+n}%
\equiv0\operatorname{mod}a_{n}$. (In the future, we shall no longer explicitly
say things like this; we shall leave it to the reader to apply Proposition
\ref{prop.mod.chain} to any chain of congruences that we write down.)

Proposition \ref{prop.mod.0} \textbf{(a)} (applied to $a_{nm+n}$ and $a_{n}$
instead of $a$ and $n$) shows that we have $a_{nm+n}\equiv0\operatorname{mod}%
a_{n}$ if and only if $a_{n}\mid a_{nm+n}$. Hence, we have $a_{n}\mid
a_{nm+n}$ (since $a_{nm+n}\equiv0\operatorname{mod}a_{n}$). In view of
$nm+n=n\left(  m+1\right)  $, this rewrites as $a_{n}\mid a_{n\left(
m+1\right)  }$. In other words, (\ref{pf.thm.rec-seq.somos-simple.b.claim})
holds for $w=m+1$. This completes the induction step. Thus,
(\ref{pf.thm.rec-seq.somos-simple.b.claim}) is proven by induction.

Therefore, Theorem \ref{thm.rec-seq.somos-simple} \textbf{(b)} is proven.

\textbf{(c)} Let $u$ and $v$ be two nonnegative integers satisfying $u\mid v$.
We must prove that $a_{u}\mid a_{v}$. If $v=0$, then this is obvious (because
if $v=0$, then $a_{v}=a_{0}=0=0a_{u}$ and therefore $a_{u}\mid a_{v}$). Hence,
for the rest of this proof, we can WLOG assume that we don't have $v=0$.
Assume this.

Thus, we don't have $v=0$. Hence, $v\neq0$, so that $v>0$ (since $v$ is nonnegative).

But $u$ divides $v$ (since $u\mid v$). In other words, there exists an integer
$w$ such that $v=uw$. Consider this $w$. If we had $w<0$, then we would have
$uw\leq0$ (since $u$ is nonnegative), which would contradict $uw=v>0$. Hence,
we cannot have $w<0$. Thus, we must have $w\geq0$. Therefore, $w\in\mathbb{N}%
$. Hence, Theorem \ref{thm.rec-seq.somos-simple} \textbf{(b)} (applied to
$n=u$) yields $a_{u}\mid a_{uw}$. In view of $v=uw$, this rewrites as
$a_{u}\mid a_{v}$. This proves Theorem \ref{thm.rec-seq.somos-simple}
\textbf{(c)}.
\end{proof}

Applying Theorem \ref{thm.rec-seq.somos-simple} \textbf{(c)} to $q=2$ and
$r=1$, we obtain the observation about divisibility made in Example
\ref{exa.rec-seq.1}.

\subsubsection{The Fibonacci sequence and a generalization}

Another example of a recursively defined sequence is the famous Fibonacci sequence:

\begin{example}
\label{exa.rec-seq.fib}The
\textit{\href{https://en.wikipedia.org/wiki/Fibonacci_number}{Fibonacci
sequence}} is the sequence $\left(  f_{0},f_{1},f_{2},\ldots\right)  $ of
integers which is defined recursively by%
\begin{align*}
f_{0}  &  =0,\ \ \ \ \ \ \ \ \ \ f_{1}=1,\ \ \ \ \ \ \ \ \ \ \text{and}\\
f_{n}  &  =f_{n-1}+f_{n-2}\ \ \ \ \ \ \ \ \ \ \text{for all }n\geq2.
\end{align*}
Let us compute its first few entries:%
\begin{align*}
f_{0}  &  =0;\\
f_{1}  &  =1;\\
f_{2}  &  =\underbrace{f_{1}}_{=1}+\underbrace{f_{0}}_{=0}=1+0=1;\\
f_{3}  &  =\underbrace{f_{2}}_{=1}+\underbrace{f_{1}}_{=1}=1+1=2;\\
f_{4}  &  =\underbrace{f_{3}}_{=2}+\underbrace{f_{2}}_{=1}=2+1=3;\\
f_{5}  &  =\underbrace{f_{4}}_{=3}+\underbrace{f_{3}}_{=2}=3+2=5;\\
f_{6}  &  =\underbrace{f_{5}}_{=5}+\underbrace{f_{4}}_{=3}=5+3=8.
\end{align*}
Again, we observe (as in Example \ref{exa.rec-seq.1}) that $f_{2}\mid f_{6}$
and $f_{3}\mid f_{6}$, which suggests that we might have $f_{u}\mid f_{v}$
whenever $u$ and $v$ are two nonnegative integers satisfying $u\mid v$.

Some further experimentation may suggest that the equality $f_{n+m+1}%
=f_{n}f_{m}+f_{n+1}f_{m+1}$ holds for all $n\in\mathbb{N}$ and $m\in
\mathbb{N}$.

Both of these conjectures will be shown in the following theorem, in greater generality.
\end{example}

\begin{theorem}
\label{thm.rec-seq.fibx}Fix some $a\in\mathbb{Z}$ and $b\in\mathbb{Z}$. Let
$\left(  x_{0},x_{1},x_{2},\ldots\right)  $ be a sequence of integers defined
recursively by%
\begin{align*}
x_{0}  &  =0,\ \ \ \ \ \ \ \ \ \ x_{1}=1,\ \ \ \ \ \ \ \ \ \ \text{and}\\
x_{n}  &  =ax_{n-1}+bx_{n-2}\ \ \ \ \ \ \ \ \ \ \text{for each }n\geq2.
\end{align*}
(Note that if $a=1$ and $b=1$, then this sequence $\left(  x_{0},x_{1}%
,x_{2},\ldots\right)  $ is precisely the Fibonacci sequence $\left(
f_{0},f_{1},f_{2},\ldots\right)  $ from Example \ref{exa.rec-seq.fib}. If
$a=0$ and $b=1$, then our sequence $\left(  x_{0},x_{1},x_{2},\ldots\right)  $
is the sequence $\left(  0,1,0,b,0,b^{2},0,b^{3},\ldots\right)  $ that
alternates between $0$'s and powers of $b$. The reader can easily work out
further examples.)

\textbf{(a)} We have $x_{n+m+1}=bx_{n}x_{m}+x_{n+1}x_{m+1}$ for all
$n\in\mathbb{N}$ and $m\in\mathbb{N}$.

\textbf{(b)} For any $n\in\mathbb{N}$ and $w\in\mathbb{N}$, we have $x_{n}\mid
x_{nw}$.

\textbf{(c)} If $u$ and $v$ are two nonnegative integers satisfying $u\mid v$,
then $x_{u}\mid x_{v}$.
\end{theorem}

Before we prove this theorem, let us discuss how \textbf{not} to prove it:

\begin{remark}
\label{rmk.ind.abstract}The proof of Theorem \ref{thm.rec-seq.fibx}
\textbf{(a)} below illustrates an important aspect of induction proofs:
Namely, when devising an induction proof, we often have not only a choice of
what variable to induct on (e.g., we could try proving Theorem
\ref{thm.rec-seq.fibx} \textbf{(a)} by induction on $n$ or by induction on
$m$), but also a choice of whether to leave the other variables fixed. For
example, let us try to prove Theorem \ref{thm.rec-seq.fibx} \textbf{(a)} by
induction on $n$ while leaving the variable $m$ fixed. That is, we fix some
$m\in\mathbb{N}$, and we define $\mathcal{A}\left(  n\right)  $ (for each
$n\in\mathbb{N}$) to be the following statement:%
\[
\left(  x_{n+m+1}=bx_{n}x_{m}+x_{n+1}x_{m+1}\right)  .
\]
Then, it is easy to check that $\mathcal{A}\left(  0\right)  $ holds, so the
induction base is complete. For the induction step, we fix some $k\in
\mathbb{N}$. (This $k$ serves the role of the \textquotedblleft$m$%
\textquotedblright\ in Theorem \ref{thm.ind.IP0}, but we cannot call it $m$
here since $m$ already stands for a fixed number.) We assume that
$\mathcal{A}\left(  k\right)  $ holds, and we intend to prove $\mathcal{A}%
\left(  k+1\right)  $.

Our induction hypothesis says that $\mathcal{A}\left(  k\right)  $ holds; in
other words, we have $x_{k+m+1}=bx_{k}x_{m}+x_{k+1}x_{m+1}$. We want to prove
$\mathcal{A}\left(  k+1\right)  $; in other words, we want to prove that
$x_{\left(  k+1\right)  +m+1}=bx_{k+1}x_{m}+x_{\left(  k+1\right)  +1}x_{m+1}$.

A short moment of deliberation shows that we cannot do this (at least not with
our current knowledge). There is no direct way of deriving $\mathcal{A}\left(
k+1\right)  $ from $\mathcal{A}\left(  k\right)  $. \textbf{However}, if we
knew that the statement $\mathcal{A}\left(  k\right)  $ holds
\textquotedblleft for $m+1$ instead of $m$\textquotedblright\ (that is, if we
knew that $x_{k+\left(  m+1\right)  +1}=bx_{k}x_{m+1}+x_{k+1}x_{\left(
m+1\right)  +1}$), then we could derive $\mathcal{A}\left(  k+1\right)  $. But
we cannot just \textquotedblleft apply $\mathcal{A}\left(  k\right)  $ to
$m+1$ instead of $m$\textquotedblright; after all, $m$ is a fixed number, so
we cannot have it take different values in $\mathcal{A}\left(  k\right)  $ and
in $\mathcal{A}\left(  k+1\right)  $.

So we are at an impasse. We got into this impasse by fixing $m$. So let us try
\textbf{not} fixing $m\in\mathbb{N}$ right away, but instead defining
$\mathcal{A}\left(  n\right)  $ (for each $n\in\mathbb{N}$) to be the
following statement:%
\[
\left(  x_{n+m+1}=bx_{n}x_{m}+x_{n+1}x_{m+1}\text{ for all }m\in
\mathbb{N}\right)  .
\]
Thus, $\mathcal{A}\left(  n\right)  $ is not a statement about a specific
integer $m$ any more, but rather a statement about all nonnegative integers
$m$. This allows us to apply $\mathcal{A}\left(  k\right)  $ to $m+1$ instead
of $m$ in the induction step. (We can still fix $m\in\mathbb{N}$
\textbf{during the induction step}; this doesn't prevent us from applying
$\mathcal{A}\left(  k\right)  $ to $m+1$ instead of $m$, since $\mathcal{A}%
\left(  k\right)  $ has been formulated before $m$ was fixed.) This way, we
arrive at the following proof:
\end{remark}

\begin{proof}
[Proof of Theorem \ref{thm.rec-seq.fibx}.]\textbf{(a)} We claim that for each
$n\in\mathbb{N}$, we have%
\begin{equation}
\left(  x_{n+m+1}=bx_{n}x_{m}+x_{n+1}x_{m+1}\text{ for all }m\in
\mathbb{N}\right)  . \label{pf.thm.rec-seq.fibx.a.claim}%
\end{equation}

Indeed, let us prove (\ref{pf.thm.rec-seq.fibx.a.claim}) by induction on $n$:

\textit{Induction base:} We have $x_{0+m+1}=bx_{0}x_{m}+x_{0+1}x_{m+1}$ for
all $m\in\mathbb{N}$\ \ \ \ \footnote{\textit{Proof.} Let $m\in\mathbb{N}$.
Then, $x_{0+m+1}=x_{m+1}$. Comparing this with $b\underbrace{x_{0}}_{=0}%
x_{m}+\underbrace{x_{0+1}}_{=x_{1}=1}x_{m+1}=b0x_{m}+1x_{m+1}=x_{m+1}$, we
obtain $x_{0+m+1}=bx_{0}x_{m}+x_{0+1}x_{m+1}$, qed.}. In other words,
(\ref{pf.thm.rec-seq.fibx.a.claim}) holds for $n=0$. This completes the
induction base.

\textit{Induction step:} Let $k\in\mathbb{N}$. Assume that
(\ref{pf.thm.rec-seq.fibx.a.claim}) holds for $n=k$. We must prove that
(\ref{pf.thm.rec-seq.fibx.a.claim}) holds for $n=k+1$.

We have assumed that (\ref{pf.thm.rec-seq.fibx.a.claim}) holds for $n=k$. In
other words, we have%
\begin{equation}
\left(  x_{k+m+1}=bx_{k}x_{m}+x_{k+1}x_{m+1}\text{ for all }m\in
\mathbb{N}\right)  . \label{pf.thm.rec-seq.fibx.a.claim.pf.IH}%
\end{equation}

Now, let $m\in\mathbb{N}$. We have $m+2\geq2$; thus, the recursive definition
of the sequence $\left(  x_{0},x_{1},x_{2},\ldots\right)  $ yields%
\begin{equation}
x_{m+2}=a\underbrace{x_{\left(  m+2\right)  -1}}_{=x_{m+1}}%
+b\underbrace{x_{\left(  m+2\right)  -2}}_{=x_{m}}=ax_{m+1}+bx_{m}.
\label{pf.thm.rec-seq.fibx.a.claim.pf.1}%
\end{equation}
The same argument (with $m$ replaced by $k$) yields%
\begin{equation}
x_{k+2}=ax_{k+1}+bx_{k}. \label{pf.thm.rec-seq.fibx.a.claim.pf.2}%
\end{equation}

But we can apply (\ref{pf.thm.rec-seq.fibx.a.claim.pf.IH}) to $m+1$ instead of
$m$. Thus, we obtain%
\begin{align*}
x_{k+\left(  m+1\right)  +1}  &  =bx_{k}x_{m+1}+x_{k+1}\underbrace{x_{\left(
m+1\right)  +1}}_{\substack{=x_{m+2}=ax_{m+1}+bx_{m}\\\text{(by
(\ref{pf.thm.rec-seq.fibx.a.claim.pf.1}))}}}\\
&  =bx_{k}x_{m+1}+\underbrace{x_{k+1}\left(  ax_{m+1}+bx_{m}\right)
}_{=ax_{k+1}x_{m+1}+bx_{k+1}x_{m}}=\underbrace{bx_{k}x_{m+1}+ax_{k+1}x_{m+1}%
}_{=\left(  ax_{k+1}+bx_{k}\right)  x_{m+1}}+bx_{k+1}x_{m}\\
&  =\underbrace{\left(  ax_{k+1}+bx_{k}\right)  }_{\substack{=x_{k+2}%
\\\text{(by (\ref{pf.thm.rec-seq.fibx.a.claim.pf.2}))}}}x_{m+1}+bx_{k+1}%
x_{m}=x_{k+2}x_{m+1}+bx_{k+1}x_{m}\\
&  =bx_{k+1}x_{m}+\underbrace{x_{k+2}}_{=x_{\left(  k+1\right)  +1}}%
x_{m+1}=bx_{k+1}x_{m}+x_{\left(  k+1\right)  +1}x_{m+1}.
\end{align*}
In view of $k+\left(  m+1\right)  +1=\left(  k+1\right)  +m+1$, this rewrites
as
\[
x_{\left(  k+1\right)  +m+1}=bx_{k+1}x_{m}+x_{\left(  k+1\right)  +1}x_{m+1}.
\]

Now, forget that we fixed $m$. We thus have shown that $x_{\left(  k+1\right)
+m+1}=bx_{k+1}x_{m}+x_{\left(  k+1\right)  +1}x_{m+1}$ for all $m\in
\mathbb{N}$. In other words, (\ref{pf.thm.rec-seq.fibx.a.claim}) holds for
$n=k+1$. This completes the induction step. Thus,
(\ref{pf.thm.rec-seq.fibx.a.claim}) is proven.

Hence, Theorem \ref{thm.rec-seq.fibx} \textbf{(a)} holds.

\textbf{(b)} Fix $n\in\mathbb{N}$. We claim that%
\begin{equation}
x_{n}\mid x_{nw}\ \ \ \ \ \ \ \ \ \ \text{for each }w\in\mathbb{N}.
\label{pf.thm.rec-seq.fibx.b.claim}%
\end{equation}

Indeed, let us prove (\ref{pf.thm.rec-seq.fibx.b.claim}) by induction on $w$:

\textit{Induction base:} We have $x_{n\cdot0}=x_{0}=0=0x_{n}$ and thus
$x_{n}\mid x_{n\cdot0}$. In other words, (\ref{pf.thm.rec-seq.fibx.b.claim})
holds for $w=0$. This completes the induction base.

\textit{Induction step:} Let $k\in\mathbb{N}$. Assume that
(\ref{pf.thm.rec-seq.fibx.b.claim}) holds for $w=k$. We must now prove that
(\ref{pf.thm.rec-seq.fibx.b.claim}) holds for $w=k+1$. In other words, we must
prove that $x_{n}\mid x_{n\left(  k+1\right)  }$.

If $n=0$, then this is true\footnote{\textit{Proof.} Let us assume that $n=0$.
Then, $x_{n\left(  k+1\right)  }=x_{0\left(  k+1\right)  }=x_{0}=0=0x_{n}$,
and thus $x_{n}\mid x_{n\left(  k+1\right)  }$, qed.}. Hence, for the rest of
this proof, we can WLOG assume that we don't have $n=0$. Assume this.

We have assumed that (\ref{pf.thm.rec-seq.fibx.b.claim}) holds for $w=k$. In
other words, we have $x_{n}\mid x_{nk}$. In other words, $x_{nk}%
\equiv0\operatorname{mod}x_{n}$.\ \ \ \ \footnote{Here, again, we have used
Proposition \ref{prop.mod.0} \textbf{(a)} (applied to $x_{nk}$ and $x_{n}$
instead of $a$ and $n$). This argument is simple enough that we will leave it
unsaid in the future.} Likewise, from $x_{n}\mid x_{n}$, we obtain
$x_{n}\equiv0\operatorname{mod}x_{n}$.

We have $n\in\mathbb{N}$ but $n\neq0$ (since we don't have $n=0$). Hence, $n$
is a positive integer. Thus, $n-1\in\mathbb{N}$. Therefore, Theorem
\ref{thm.rec-seq.fibx} \textbf{(a)} (applied to $nk$ and $n-1$ instead of $n$
and $m$) yields
\[
x_{nk+\left(  n-1\right)  +1}=bx_{nk}x_{n-1}+x_{nk+1}x_{\left(  n-1\right)
+1}.
\]
In view of $nk+\left(  n-1\right)  +1=n\left(  k+1\right)  $, this rewrites as%
\[
x_{n\left(  k+1\right)  }=b\underbrace{x_{nk}}_{\equiv0\operatorname{mod}%
x_{n}}x_{n-1}+x_{nk+1}\underbrace{x_{\left(  n-1\right)  +1}}_{=x_{n}%
\equiv0\operatorname{mod}x_{n}}\equiv b0x_{n-1}+x_{nk+1}0=0\operatorname{mod}%
x_{n}.
\]
\footnote{We have used substitutivity for congruences in this computation.
Here is, again, a way to rewrite it without this use:
\par
We have $x_{n\left(  k+1\right)  }=bx_{nk}x_{n-1}+x_{nk+1}x_{\left(
n-1\right)  +1}$. But $b\equiv b\operatorname{mod}x_{n}$ (by Proposition
\ref{prop.mod.transi} \textbf{(a)}) and $x_{n-1}\equiv x_{n-1}%
\operatorname{mod}x_{n}$ (for the same reason) and $x_{nk+1}\equiv
x_{nk+1}\operatorname{mod}x_{n}$ (for the same reason). Now, Proposition
\ref{prop.mod.+-*} \textbf{(c)} (applied to $b$, $b$, $x_{nk}$, $0$ and
$x_{n}$ instead of $a$, $b$, $c$, $d$ and $n$) yields $bx_{nk}\equiv
b0\operatorname{mod}x_{n}$ (since $b\equiv b\operatorname{mod}x_{n}$ and
$x_{nk}\equiv0\operatorname{mod}x_{n}$). Hence, Proposition \ref{prop.mod.+-*}
\textbf{(c)} (applied to $bx_{nk}$, $b0$, $x_{n-1}$, $x_{n-1}$ and $x_{n}$
instead of $a$, $b$, $c$, $d$ and $n$) yields $bx_{nk}x_{n-1}\equiv
b0x_{n-1}\operatorname{mod}x_{n}$ (since $bx_{nk}\equiv b0\operatorname{mod}%
x_{n}$ and $x_{n-1}\equiv x_{n-1}\operatorname{mod}x_{n}$). Also,
$x_{nk+1}x_{\left(  n-1\right)  +1}\equiv x_{nk+1}x_{\left(  n-1\right)
+1}\operatorname{mod}x_{n}$ (by Proposition \ref{prop.mod.transi}
\textbf{(a)}). Hence, Proposition \ref{prop.mod.+-*} \textbf{(a)} (applied to
$bx_{nk}x_{n-1}$, $b0x_{n-1}$, $x_{nk+1}x_{\left(  n-1\right)  +1}$,
$x_{nk+1}x_{\left(  n-1\right)  +1}$ and $x_{n}$ instead of $a$, $b$, $c$, $d$
and $n$) yields
\[
bx_{nk}x_{n-1}+x_{nk+1}x_{\left(  n-1\right)  +1}\equiv b0x_{n-1}%
+x_{nk+1}x_{\left(  n-1\right)  +1}\operatorname{mod}x_{n}%
\]
(since $bx_{nk}x_{n-1}\equiv b0x_{n-1}\operatorname{mod}x_{n}$ and
$x_{nk+1}x_{\left(  n-1\right)  +1}\equiv x_{nk+1}x_{\left(  n-1\right)
+1}\operatorname{mod}x_{n}$).
\par
Also, Proposition \ref{prop.mod.+-*} \textbf{(c)} (applied to $x_{nk+1}$,
$x_{nk+1}$, $x_{\left(  n-1\right)  +1}$, $0$ and $x_{n}$ instead of $a$, $b$,
$c$, $d$ and $n$) yields $x_{nk+1}x_{\left(  n-1\right)  +1}\equiv
x_{nk+1}0\operatorname{mod}x_{n}$ (since $x_{nk+1}\equiv x_{nk+1}%
\operatorname{mod}x_{n}$ and $x_{\left(  n-1\right)  +1}=x_{n}\equiv
0\operatorname{mod}x_{n}$). Furthermore, $b0x_{n-1}\equiv b0x_{n-1}%
\operatorname{mod}x_{n}$ (by Proposition \ref{prop.mod.transi} \textbf{(a)}).
Finally, Proposition \ref{prop.mod.+-*} \textbf{(a)} (applied to $b0x_{n-1}$,
$b0x_{n-1}$, $x_{nk+1}x_{\left(  n-1\right)  +1}$, $x_{nk+1}0$ and $x_{n}$
instead of $a$, $b$, $c$, $d$ and $n$) yields
\[
b0x_{n-1}+x_{nk+1}x_{\left(  n-1\right)  +1}\equiv b0x_{n-1}+x_{nk+1}%
0\operatorname{mod}x_{n}%
\]
(since $b0x_{n-1}\equiv b0x_{n-1}\operatorname{mod}x_{n}$ and $x_{nk+1}%
x_{\left(  n-1\right)  +1}\equiv x_{nk+1}0\operatorname{mod}x_{n}$). Thus,%
\begin{align*}
x_{n\left(  k+1\right)  }  &  =bx_{nk}x_{n-1}+x_{nk+1}x_{\left(  n-1\right)
+1}\equiv b0x_{n-1}+x_{nk+1}x_{\left(  n-1\right)  +1}\\
&  \equiv b0x_{n-1}+x_{nk+1}0=0\operatorname{mod}x_{n}.
\end{align*}
So we have proven that $x_{n\left(  k+1\right)  }\equiv0\operatorname{mod}%
x_{n}$.} Thus, we have shown that $x_{n\left(  k+1\right)  }\equiv
0\operatorname{mod}x_{n}$. In other words, $x_{n}\mid x_{n\left(  k+1\right)
}$ (again, this follows from Proposition \ref{prop.mod.0} \textbf{(a)}). In
other words, (\ref{pf.thm.rec-seq.fibx.b.claim}) holds for $w=k+1$. This
completes the induction step. Hence, (\ref{pf.thm.rec-seq.fibx.b.claim}) is
proven by induction.

This proves Theorem \ref{thm.rec-seq.fibx} \textbf{(b)}.

\begin{vershort}
\textbf{(c)} Theorem \ref{thm.rec-seq.fibx} \textbf{(c)} can be derived from
Theorem \ref{thm.rec-seq.fibx} \textbf{(b)} in the same way as Theorem
\ref{thm.rec-seq.somos-simple} \textbf{(c)} was derived from Theorem
\ref{thm.rec-seq.somos-simple} \textbf{(b)}. \qedhere

\end{vershort}

\begin{verlong}
\textbf{(c)} Let $u$ and $v$ be two nonnegative integers satisfying $u\mid v$.
We must prove that $x_{u}\mid x_{v}$. If $v=0$, then this is obvious (because
if $v=0$, then $x_{v}=x_{0}=0=0x_{u}$ and therefore $x_{u}\mid x_{v}$). Hence,
for the rest of this proof, we can WLOG assume that we don't have $v=0$.
Assume this.

Thus, we don't have $v=0$. Hence, $v\neq0$, so that $v>0$ (since $v$ is nonnegative).

But $u$ divides $v$ (since $u\mid v$). In other words, there exists an integer
$w$ such that $v=uw$. Consider this $w$. If we had $w<0$, then we would have
$uw\leq0$ (since $u$ is nonnegative), which would contradict $uw=v>0$. Hence,
we cannot have $w<0$. Thus, we must have $w\geq0$. Therefore, $w\in\mathbb{N}%
$. Hence, Theorem \ref{thm.rec-seq.fibx} \textbf{(b)} (applied to $n=u$)
yields $x_{u}\mid x_{uw}$. In view of $v=uw$, this rewrites as $x_{u}\mid
x_{v}$. This proves Theorem \ref{thm.rec-seq.fibx} \textbf{(c)}.
\end{verlong}
\end{proof}

Applying Theorem \ref{thm.rec-seq.fibx} \textbf{(a)} to $a=1$ and $b=1$, we
obtain the equality $f_{n+m+1}=f_{n}f_{m}+f_{n+1}f_{m+1}$ noticed in Example
\ref{exa.rec-seq.fib}. Applying Theorem \ref{thm.rec-seq.fibx} \textbf{(c)} to
$a=1$ and $b=1$, we obtain the observation about divisibility made in Example
\ref{exa.rec-seq.fib}.

Note that part \textbf{(a)} of Theorem \ref{thm.rec-seq.fibx} still works if
$a$ and $b$ are real numbers (instead of being integers). But of course, in
this case, $\left(  x_{0},x_{1},x_{2},\ldots\right)  $ will be merely a
sequence of real numbers (rather than a sequence of integers), and thus parts
\textbf{(b)} and \textbf{(c)} of Theorem \ref{thm.rec-seq.fibx} will no longer
make sense (since divisibility is only defined for integers).

\subsection{\label{sect.ind.trinum}The sum of the first $n$ positive integers}

We now come to one of the most classical examples of a proof by induction:
Namely, we shall prove the fact that for each $n\in\mathbb{N}$, the sum of the
first $n$ positive integers (that is, the sum $1+2+\cdots+n$) equals
$\dfrac{n\left(  n+1\right)  }{2}$. However, there is a catch here, which is
easy to overlook if one isn't trying to be completely rigorous: We don't
really know yet whether there is such a thing as \textquotedblleft the sum of
the first $n$ positive integers\textquotedblright! To be more precise, we have
introduced the $\sum$ sign in Section \ref{sect.sums-repetitorium}, which
would allow us to define the sum of the first $n$ positive integers (as
$\sum_{i=1}^{n}i$); but our definition of the $\sum$ sign relied on a fact
which we have not proved yet (namely, the fact that the right hand side of
(\ref{eq.sum.def.1}) does not depend on the choice of $t$). We shall prove
this fact later (Theorem \ref{thm.ind.gen-com.wd} \textbf{(a)}), but for now
we prefer not to use it. Instead, let us replace the notion of
\textquotedblleft the sum of the first $n$ positive integers\textquotedblright%
\ by a recursively defined sequence:

\begin{proposition}
\label{prop.rec-seq.triangular}Let $\left(  t_{0},t_{1},t_{2},\ldots\right)  $
be a sequence of integers defined recursively by%
\begin{align*}
t_{0}  &  =0,\ \ \ \ \ \ \ \ \ \ \text{and}\\
t_{n}  &  =t_{n-1}+n\ \ \ \ \ \ \ \ \ \ \text{for each }n\geq1.
\end{align*}
Then,%
\begin{equation}
t_{n}=\dfrac{n\left(  n+1\right)  }{2}\ \ \ \ \ \ \ \ \ \ \text{for each }%
n\in\mathbb{N}. \label{eq.prop.rec-seq.triangular.claim}%
\end{equation}

\end{proposition}

The sequence $\left(  t_{0},t_{1},t_{2},\ldots\right)  $ defined in
Proposition \ref{prop.rec-seq.triangular} is known as the \textit{sequence of
triangular numbers}. Its definition shows that%
\begin{align*}
t_{0}  &  =0;\\
t_{1}  &  =\underbrace{t_{0}}_{=0}+1=0+1=1;\\
t_{2}  &  =\underbrace{t_{1}}_{=1}+2=1+2;\\
t_{3}  &  =\underbrace{t_{2}}_{=1+2}+3=\left(  1+2\right)  +3;\\
t_{4}  &  =\underbrace{t_{3}}_{=\left(  1+2\right)  +3}+4=\left(  \left(
1+2\right)  +3\right)  +4;\\
t_{5}  &  =\underbrace{t_{4}}_{=\left(  \left(  1+2\right)  +3\right)
+4}+5=\left(  \left(  \left(  1+2\right)  +3\right)  +4\right)  +5
\end{align*}
\footnote{Note that we write \textquotedblleft$\left(  \left(  \left(
1+2\right)  +3\right)  +4\right)  +5$\textquotedblright\ and not
\textquotedblleft$1+2+3+4+5$\textquotedblright. The reason for this is that we
haven't proven yet that the expression \textquotedblleft$1+2+3+4+5$%
\textquotedblright\ is well-defined. (This expression \textbf{is}
well-defined, but this will only be clear once we have proven Theorem
\ref{thm.ind.gen-com.wd} \textbf{(a)} below.)} and so on; this explains why it
makes sense to think of $t_{n}$ as the sum of the first $n$ positive integers.
(This is legitimate even when $n=0$, because the sum of the first $0$ positive
integers is an empty sum, and an empty sum is always defined to be equal to
$0$.) Once we have convinced ourselves that \textquotedblleft the sum of the
first $n$ positive integers\textquotedblright\ is a well-defined concept, it
will be easy to see (by induction) that $t_{n}$ \textbf{is} the sum of the
first $n$ positive integers whenever $n\in\mathbb{N}$. Therefore, Proposition
\ref{prop.rec-seq.triangular} will tell us that the sum of the first $n$
positive integers equals $\dfrac{n\left(  n+1\right)  }{2}$ whenever
$n\in\mathbb{N}$.

For now, let us prove Proposition \ref{prop.rec-seq.triangular}:

\begin{proof}
[Proof of Proposition \ref{prop.rec-seq.triangular}.]We shall prove
(\ref{eq.prop.rec-seq.triangular.claim}) by induction on $n$:

\textit{Induction base:} Comparing $t_{0}=0$ with $\dfrac{0\left(  0+1\right)
}{2}=0$, we obtain $t_{0}=\dfrac{0\left(  0+1\right)  }{2}$. In other words,
(\ref{eq.prop.rec-seq.triangular.claim}) holds for $n=0$. This completes the
induction base.

\textit{Induction step:} Let $m\in\mathbb{N}$. Assume that
(\ref{eq.prop.rec-seq.triangular.claim}) holds for $n=m$. We must prove that
(\ref{eq.prop.rec-seq.triangular.claim}) holds for $n=m+1$.

We have assumed that (\ref{eq.prop.rec-seq.triangular.claim}) holds for $n=m$.
In other words, we have $t_{m}=\dfrac{m\left(  m+1\right)  }{2}$.

Recall that $t_{n}=t_{n-1}+n$ for each $n\geq1$. Applying this to $n=m+1$, we
obtain%
\begin{align*}
t_{m+1}  &  =\underbrace{t_{\left(  m+1\right)  -1}}_{=t_{m}=\dfrac{m\left(
m+1\right)  }{2}}+\left(  m+1\right)  =\dfrac{m\left(  m+1\right)  }%
{2}+\left(  m+1\right)  =\dfrac{m\left(  m+1\right)  +2\left(  m+1\right)
}{2}\\
&  =\dfrac{\left(  m+2\right)  \left(  m+1\right)  }{2}=\dfrac{\left(
m+1\right)  \left(  m+2\right)  }{2}=\dfrac{\left(  m+1\right)  \left(
\left(  m+1\right)  +1\right)  }{2}%
\end{align*}
(since $m+2=\left(  m+1\right)  +1$). In other words,
(\ref{eq.prop.rec-seq.triangular.claim}) holds for $n=m+1$. This completes the
induction step. Hence, (\ref{eq.prop.rec-seq.triangular.claim}) is proven by
induction. This proves Proposition \ref{prop.rec-seq.triangular}.
\end{proof}

\subsection{\label{sect.ind.max}Induction on a derived quantity: maxima of
sets}

\subsubsection{Defining maxima}

We have so far been applying the Induction Principle in fairly obvious ways:
With the exception of our proof of Proposition \ref{prop.mod.chain}, we have
mostly been doing induction on a variable ($n$ or $k$ or $i$) that already
appeared in the claim that we were proving. But sometimes, it is worth doing
induction on a variable that does \textbf{not} explicitly appear in this claim
(which, formally speaking, means that we introduce a new variable to do
induction on). For example, the claim might be saying \textquotedblleft Each
nonempty finite set $S$ of integers has a largest element\textquotedblright,
and we prove it by induction on $\left\vert S\right\vert -1$. This means that
instead of directly proving the claim itself, we rather prove the equivalent
claim \textquotedblleft For each $n\in\mathbb{N}$, each nonempty finite set
$S$ of integers satisfying $\left\vert S\right\vert -1=n$ has a largest
element\textquotedblright\ by induction on $n$. We shall show this proof in
more detail below (see Theorem \ref{thm.ind.max}). First, we prepare by
discussing largest elements of sets in general.

\begin{definition}
Let $S$ be a set of integers (or rational numbers, or real numbers). A
\textit{maximum} of $S$ is defined to be an element $s\in S$ that satisfies%
\[
\left(  s\geq t\text{ for each }t\in S\right)  .
\]
In other words, a maximum of $S$ is defined to be an element of $S$ which is
greater or equal to each element of $S$.

(The plural of the word \textquotedblleft maximum\textquotedblright\ is
\textquotedblleft maxima\textquotedblright.)
\end{definition}

\begin{example}
The set $\left\{  2,4,5\right\}  $ has exactly one maximum: namely, $5$.

The set $\mathbb{N}=\left\{  0,1,2,\ldots\right\}  $ has no maximum: If $k$
was a maximum of $\mathbb{N}$, then we would have $k\geq k+1$, which is absurd.

The set $\left\{  0,-1,-2,\ldots\right\}  $ has a maximum: namely, $0$.

The set $\varnothing$ has no maximum, since a maximum would have to be an
element of $\varnothing$.
\end{example}

In Theorem \ref{thm.ind.max}, we shall soon show that every nonempty finite
set of integers has a maximum. First, we prove that a maximum is unique if it exists:

\begin{proposition}
\label{prop.ind.max-uni}Let $S$ be a set of integers (or rational numbers, or
real numbers). Then, $S$ has \textbf{at most one} maximum.
\end{proposition}

\begin{proof}
[Proof of Proposition \ref{prop.ind.max-uni}.]Let $s_{1}$ and $s_{2}$ be two
maxima of $S$. We shall show that $s_{1}=s_{2}$.

Indeed, $s_{1}$ is a maximum of $S$. In other words, $s_{1}$ is an element
$s\in S$ that satisfies $\left(  s\geq t\text{ for each }t\in S\right)  $ (by
the definition of a maximum). In other words, $s_{1}$ is an element of $S$ and
satisfies%
\begin{equation}
\left(  s_{1}\geq t\text{ for each }t\in S\right)  .
\label{pf.prop.ind.max-uni.1}%
\end{equation}
The same argument (applied to $s_{2}$ instead of $s_{1}$) shows that $s_{2}$
is an element of $S$ and satisfies%
\begin{equation}
\left(  s_{2}\geq t\text{ for each }t\in S\right)  .
\label{pf.prop.ind.max-uni.2}%
\end{equation}

Now, $s_{1}$ is an element of $S$. Hence, (\ref{pf.prop.ind.max-uni.2})
(applied to $t=s_{1}$) yields $s_{2}\geq s_{1}$. But the same argument (with
the roles of $s_{1}$ and $s_{2}$ interchanged) shows that $s_{1}\geq s_{2}$.
Combining this with $s_{2}\geq s_{1}$, we obtain $s_{1}=s_{2}$.

Now, forget that we fixed $s_{1}$ and $s_{2}$. We thus have shown that if
$s_{1}$ and $s_{2}$ are two maxima of $S$, then $s_{1}=s_{2}$. In other words,
any two maxima of $S$ are equal. In other words, $S$ has \textbf{at most one}
maximum. This proves Proposition \ref{prop.ind.max-uni}.
\end{proof}

\begin{definition}
Let $S$ be a set of integers (or rational numbers, or real numbers).
Proposition \ref{prop.ind.max-uni} shows that $S$ has \textbf{at most one}
maximum. Thus, if $S$ has a maximum, then this maximum is the unique maximum
of $S$; we shall thus call it \textit{the maximum} of $S$ or \textit{the
largest element} of $S$. We shall denote this maximum by $\max S$.
\end{definition}

Thus, if $S$ is a set of integers (or rational numbers, or real numbers) that
has a maximum, then this maximum $\max S$ satisfies
\begin{equation}
\max S\in S \label{eq.ind.max.def-max.1}%
\end{equation}
and%
\begin{equation}
\left(  \max S\geq t\text{ for each }t\in S\right)
\label{eq.ind.max.def-max.2}%
\end{equation}
(because of the definition of a maximum).

Let us next show two simple facts:

\begin{lemma}
\label{lem.ind.max-1el}Let $x$ be an integer (or rational number, or real
number). Then, the set $\left\{  x\right\}  $ has a maximum, namely $x$.
\end{lemma}

\begin{proof}
[Proof of Lemma \ref{lem.ind.max-1el}.]Clearly, $x\geq x$. Thus, $x\geq t$ for
each $t\in\left\{  x\right\}  $ (because the only $t\in\left\{  x\right\}  $
is $x$). In other words, $x$ is an element $s\in\left\{  x\right\}  $ that
satisfies \newline$\left(  s\geq t\text{ for each }t\in\left\{  x\right\}
\right)  $ (since $x\in\left\{  x\right\}  $).

But recall that a maximum of $\left\{  x\right\}  $ means an element
$s\in\left\{  x\right\}  $ that satisfies $\left(  s\geq t\text{ for each
}t\in\left\{  x\right\}  \right)  $ (by the definition of a maximum). Hence,
$x$ is a maximum of $\left\{  x\right\}  $ (since $x$ is such an element).
Thus, the set $\left\{  x\right\}  $ has a maximum, namely $x$. This proves
Lemma \ref{lem.ind.max-1el}.
\end{proof}

\begin{proposition}
\label{prop.ind.max-PuQ}Let $P$ and $Q$ be two sets of integers (or rational
numbers, or real numbers). Assume that $P$ has a maximum, and assume that $Q$
has a maximum. Then, the set $P\cup Q$ has a maximum.
\end{proposition}

\begin{proof}
[Proof of Proposition \ref{prop.ind.max-PuQ}.]We know that $P$ has a maximum;
it is denoted by $\max P$. We also know that $Q$ has a maximum; it is denoted
by $\max Q$. The sets $P$ and $Q$ play symmetric roles in Proposition
\ref{prop.ind.max-PuQ} (since $P\cup Q=Q\cup P$). Thus, we can WLOG assume
that $\max P\geq\max Q$ (since otherwise, we can simply swap $P$ with $Q$,
without altering the meaning of Proposition \ref{prop.ind.max-PuQ}). Assume this.

Now, (\ref{eq.ind.max.def-max.1}) (applied to $S=P$) shows that $\max P\in
P\subseteq P\cup Q$. Furthermore, we claim that%
\begin{equation}
\left(  \max P\geq t\text{ for each }t\in P\cup Q\right)  .
\label{pf.prop.ind.max-PuQ.1}%
\end{equation}

[\textit{Proof of (\ref{pf.prop.ind.max-PuQ.1}):} Let $t\in P\cup Q$. We must
show that $\max P\geq t$.

We have $t\in P\cup Q$. In other words, $t\in P$ or $t\in Q$. Hence, we are in
one of the following two cases:

\textit{Case 1:} We have $t\in P$.

\textit{Case 2:} We have $t\in Q$.

(These two cases might have overlap, but there is nothing wrong about this.)

Let us first consider Case 1. In this case, we have $t\in P$. Hence,
(\ref{eq.ind.max.def-max.2}) (applied to $S=P$) yields $\max P\geq t$. Hence,
$\max P\geq t$ is proven in Case 1.

Let us next consider Case 2. In this case, we have $t\in Q$. Hence,
(\ref{eq.ind.max.def-max.2}) (applied to $S=Q$) yields $\max Q\geq t$. Hence,
$\max P\geq\max Q\geq t$. Thus, $\max P\geq t$ is proven in Case 2.

We have now proven $\max P\geq t$ in each of the two Cases 1 and 2. Since
these two Cases cover all possibilities, we thus conclude that $\max P\geq t$
always holds. This proves (\ref{pf.prop.ind.max-PuQ.1}).]

Now, $\max P$ is an element $s\in P\cup Q$ that satisfies $\left(  s\geq
t\text{ for each }t\in P\cup Q\right)  $ (since $\max P\in P\cup Q$ and
$\left(  \max P\geq t\text{ for each }t\in P\cup Q\right)  $).

But recall that a maximum of $P\cup Q$ means an element $s\in P\cup Q$ that
satisfies $\left(  s\geq t\text{ for each }t\in P\cup Q\right)  $ (by the
definition of a maximum). Hence, $\max P$ is a maximum of $P\cup Q$ (since
$\max P$ is such an element). Thus, the set $P\cup Q$ has a maximum. This
proves Proposition \ref{prop.ind.max-PuQ}.
\end{proof}

\subsubsection{Nonempty finite sets of integers have maxima}

\begin{theorem}
\label{thm.ind.max}Let $S$ be a nonempty finite set of integers. Then, $S$ has
a maximum.
\end{theorem}

\begin{proof}
[First proof of Theorem \ref{thm.ind.max}.]First of all, let us forget that we
fixed $S$. So we want to prove that if $S$ is a nonempty finite set of
integers, then $S$ has a maximum.

For each $n\in\mathbb{N}$, we let $\mathcal{A}\left(  n\right)  $ be the
statement%
\[
\left(
\begin{array}
[c]{c}%
\text{if }S\text{ is a nonempty finite set of integers satisfying }\left\vert
S\right\vert -1=n\text{,}\\
\text{then }S\text{ has a maximum}%
\end{array}
\right)  .
\]

We claim that $\mathcal{A}\left(  n\right)  $ holds for all $n\in\mathbb{N}$.

Indeed, let us prove this by induction on $n$:

\textit{Induction base:} If $S$ is a nonempty finite set of integers
satisfying $\left\vert S\right\vert -1=0$, then $S$ has a
maximum\footnote{\textit{Proof.} Let $S$ be a nonempty finite set of integers
satisfying $\left\vert S\right\vert -1=0$. We must show that $S$ has a
maximum.
\par
Indeed, $\left\vert S\right\vert =1$ (since $\left\vert S\right\vert -1=0$).
In other words, $S$ is a $1$-element set. In other words, $S=\left\{
x\right\}  $ for some integer $x$. Consider this $x$. Lemma
\ref{lem.ind.max-1el} shows that the set $\left\{  x\right\}  $ has a maximum.
In other words, the set $S$ has a maximum (since $S=\left\{  x\right\}  $).
This completes our proof.}. But this is exactly the statement $\mathcal{A}%
\left(  0\right)  $. Hence, $\mathcal{A}\left(  0\right)  $ holds. This
completes the induction base.

\textit{Induction step:} Let $m\in\mathbb{N}$. Assume that $\mathcal{A}\left(
m\right)  $ holds. We shall now show that $\mathcal{A}\left(  m+1\right)  $ holds.

We have assumed that $\mathcal{A}\left(  m\right)  $ holds. In other words,%
\begin{equation}
\left(
\begin{array}
[c]{c}%
\text{if }S\text{ is a nonempty finite set of integers satisfying }\left\vert
S\right\vert -1=m\text{,}\\
\text{then }S\text{ has a maximum}%
\end{array}
\right)  \label{pf.thm.ind.max.IH}%
\end{equation}
(because this is what the statement $\mathcal{A}\left(  m\right)  $ says).

Now, let $S$ be a nonempty finite set of integers satisfying $\left\vert
S\right\vert -1=m+1$. There exists some $t\in S$ (since $S$ is nonempty).
Consider this $t$. We have $\left(  S\setminus\left\{  t\right\}  \right)
\cup\left\{  t\right\}  =S\cup\left\{  t\right\}  =S$ (since $t\in S$).

From $t\in S$, we obtain $\left\vert S\setminus\left\{  t\right\}  \right\vert
=\left\vert S\right\vert -1=m+1>m\geq0$ (since $m\in\mathbb{N}$). Hence, the
set $S\setminus\left\{  t\right\}  $ is nonempty. Furthermore, this set
$S\setminus\left\{  t\right\}  $ is finite (since $S$ is finite) and satisfies
$\left\vert S\setminus\left\{  t\right\}  \right\vert -1=m$ (since $\left\vert
S\setminus\left\{  t\right\}  \right\vert =m+1$). Hence,
(\ref{pf.thm.ind.max.IH}) (applied to $S\setminus\left\{  t\right\}  $ instead
of $S$) shows that $S\setminus\left\{  t\right\}  $ has a maximum. Also, Lemma
\ref{lem.ind.max-1el} (applied to $x=t$) shows that the set $\left\{
t\right\}  $ has a maximum, namely $t$. Hence, Proposition
\ref{prop.ind.max-PuQ} (applied to $P=S\setminus\left\{  t\right\}  $ and
$Q=\left\{  t\right\}  $) shows that the set $\left(  S\setminus\left\{
t\right\}  \right)  \cup\left\{  t\right\}  $ has a maximum. Since $\left(
S\setminus\left\{  t\right\}  \right)  \cup\left\{  t\right\}  =S$, this
rewrites as follows: The set $S$ has a maximum.

Now, forget that we fixed $S$. We thus have shown that if $S$ is a nonempty
finite set of integers satisfying $\left\vert S\right\vert -1=m+1$, then $S$
has a maximum. But this is precisely the statement $\mathcal{A}\left(
m+1\right)  $. Hence, we have shown that $\mathcal{A}\left(  m+1\right)  $
holds. This completes the induction step.

Thus, we have proven (by induction) that $\mathcal{A}\left(  n\right)  $ holds
for all $n\in\mathbb{N}$. In other words, for all $n\in\mathbb{N}$, the
following holds:%
\begin{equation}
\left(
\begin{array}
[c]{c}%
\text{if }S\text{ is a nonempty finite set of integers satisfying }\left\vert
S\right\vert -1=n\text{,}\\
\text{then }S\text{ has a maximum}%
\end{array}
\right)  \label{pf.thm.ind.max.AT}%
\end{equation}
(because this is what $\mathcal{A}\left(  n\right)  $ says).

Now, let $S$ be a nonempty finite set of integers. We shall prove that $S$ has
a maximum.

Indeed, $\left\vert S\right\vert \in\mathbb{N}$ (since $S$ is finite) and
$\left\vert S\right\vert >0$ (since $S$ is nonempty); hence, $\left\vert
S\right\vert \geq1$. Thus, $\left\vert S\right\vert -1\geq0$, so that
$\left\vert S\right\vert -1\in\mathbb{N}$. Hence, we can define $n\in
\mathbb{N}$ by $n=\left\vert S\right\vert -1$. Consider this $n$. Thus,
$\left\vert S\right\vert -1=n$. Hence, (\ref{pf.thm.ind.max.AT}) shows that
$S$ has a maximum. This proves Theorem \ref{thm.ind.max}.
\end{proof}

\subsubsection{Conventions for writing induction proofs on derived quantities}

Let us take a closer look at the proof we just gave. The definition of the
statement $\mathcal{A}\left(  n\right)  $ was not exactly unmotivated: This
statement simply says that Theorem \ref{thm.ind.max} holds under the condition
that $\left\vert S\right\vert -1=n$. Thus, by introducing $\mathcal{A}\left(
n\right)  $, we have \textquotedblleft sliced\textquotedblright\ Theorem
\ref{thm.ind.max} into a sequence of statements $\mathcal{A}\left(  0\right)
,\mathcal{A}\left(  1\right)  ,\mathcal{A}\left(  2\right)  ,\ldots$, which
then allowed us to prove these statements by induction on $n$ even though no
\textquotedblleft$n$\textquotedblright\ appeared in Theorem \ref{thm.ind.max}
itself. This kind of strategy applies to various other problems. Again, we
don't need to explicitly define the statement $\mathcal{A}\left(  n\right)  $
if it is simply saying that the claim we are trying to prove (in our case,
Theorem \ref{thm.ind.max}) holds under the condition that $\left\vert
S\right\vert -1=n$; we can just say that we are doing \textquotedblleft
induction on $\left\vert S\right\vert -1$\textquotedblright. More generally:

\begin{convention}
\label{conv.ind.IP0der}Let $\mathcal{B}$ be a logical statement that involves
some variables $v_{1},v_{2},v_{3},\ldots$. (For example, $\mathcal{B}$ can be
the statement of Theorem \ref{thm.ind.max}; then, there is only one variable,
namely $S$.)

Let $q$ be some expression (involving the variables $v_{1},v_{2},v_{3},\ldots$
or some of them) that has the property that whenever the variables
$v_{1},v_{2},v_{3},\ldots$ satisfy the assumptions of $\mathcal{B}$, the
expression $q$ evaluates to some nonnegative integer. (For example, if
$\mathcal{B}$ is the statement of Theorem \ref{thm.ind.max}, then $q$ can be
the expression $\left\vert S\right\vert -1$, because it is easily seen that if
$S$ is a nonempty finite set of integers, then $\left\vert S\right\vert -1$ is
a nonnegative integer.)

Assume that you want to prove the statement $\mathcal{B}$. Then, you can
proceed as follows: For each $n\in\mathbb{N}$, define $\mathcal{A}\left(
n\right)  $ to be the statement saying that\footnotemark%
\[
\left(  \text{the statement }\mathcal{B}\text{ holds under the condition that
}q=n\right)  .
\]
Then, prove $\mathcal{A}\left(  n\right)  $ by induction on $n$. Thus:

\begin{itemize}
\item The \textit{induction base} consists in proving that the statement
$\mathcal{B}$ holds under the condition that $q=0$.

\item The \textit{induction step} consists in fixing $m\in\mathbb{N}$, and
showing that if the statement $\mathcal{B}$ holds under the condition that
$q=m$, then the statement $\mathcal{B}$ holds under the condition that $q=m+1$.
\end{itemize}

Once this induction proof is finished, it immediately follows that the
statement $\mathcal{B}$ always holds.

This strategy of proof is called \textquotedblleft induction on $q$%
\textquotedblright\ (or \textquotedblleft induction over $q$\textquotedblright%
). Once you have specified what $q$ is, you don't need to explicitly define
$\mathcal{A}\left(  n\right)  $, nor do you ever need to mention $n$.
\end{convention}

\footnotetext{We assume that no variable named \textquotedblleft%
$n$\textquotedblright\ appears in the statement $\mathcal{B}$; otherwise, we
need a different letter for our new variable in order to avoid confusion.}%
Using this convention, we can rewrite our above proof of Theorem
\ref{thm.ind.max} as follows:

\begin{proof}
[First proof of Theorem \ref{thm.ind.max} (second version).]It is easy to see
that $\left\vert S\right\vert -1\in\mathbb{N}$%
\ \ \ \ \footnote{\textit{Proof.} We have $\left\vert S\right\vert
\in\mathbb{N}$ (since $S$ is finite) and $\left\vert S\right\vert >0$ (since
$S$ is nonempty); hence, $\left\vert S\right\vert \geq1$. Thus, $\left\vert
S\right\vert -1\in\mathbb{N}$, qed.}. Hence, we can apply induction on
$\left\vert S\right\vert -1$ to prove Theorem \ref{thm.ind.max}:

\textit{Induction base:} Theorem \ref{thm.ind.max} holds under the condition
that $\left\vert S\right\vert -1=0$\ \ \ \ \footnote{\textit{Proof.} Let $S$
be as in Theorem \ref{thm.ind.max}, and assume that $\left\vert S\right\vert
-1=0$. We must show that the claim of Theorem \ref{thm.ind.max} holds.
\par
Indeed, $\left\vert S\right\vert =1$ (since $\left\vert S\right\vert -1=0$).
In other words, $S$ is a $1$-element set. In other words, $S=\left\{
x\right\}  $ for some integer $x$. Consider this $x$. Lemma
\ref{lem.ind.max-1el} shows that the set $\left\{  x\right\}  $ has a maximum.
In other words, the set $S$ has a maximum (since $S=\left\{  x\right\}  $). In
other words, the claim of Theorem \ref{thm.ind.max} holds. This completes our
proof.}. This completes the induction base.

\textit{Induction step:} Let $m\in\mathbb{N}$. Assume that Theorem
\ref{thm.ind.max} holds under the condition that $\left\vert S\right\vert
-1=m$. We shall now show that Theorem \ref{thm.ind.max} holds under the
condition that $\left\vert S\right\vert -1=m+1$.

We have assumed that Theorem \ref{thm.ind.max} holds under the condition that
$\left\vert S\right\vert -1=m$. In other words,%
\begin{equation}
\left(
\begin{array}
[c]{c}%
\text{if }S\text{ is a nonempty finite set of integers satisfying }\left\vert
S\right\vert -1=m\text{,}\\
\text{then }S\text{ has a maximum}%
\end{array}
\right)  . \label{pf.thm.ind.max.ver2.1}%
\end{equation}

Now, let $S$ be a nonempty finite set of integers satisfying $\left\vert
S\right\vert -1=m+1$. There exists some $t\in S$ (since $S$ is nonempty).
Consider this $t$. We have $\left(  S\setminus\left\{  t\right\}  \right)
\cup\left\{  t\right\}  =S\cup\left\{  t\right\}  =S$ (since $t\in S$).

From $t\in S$, we obtain $\left\vert S\setminus\left\{  t\right\}  \right\vert
=\left\vert S\right\vert -1=m+1>m\geq0$ (since $m\in\mathbb{N}$). Hence, the
set $S\setminus\left\{  t\right\}  $ is nonempty. Furthermore, this set
$S\setminus\left\{  t\right\}  $ is finite (since $S$ is finite) and satisfies
$\left\vert S\setminus\left\{  t\right\}  \right\vert -1=m$ (since $\left\vert
S\setminus\left\{  t\right\}  \right\vert =m+1$). Hence,
(\ref{pf.thm.ind.max.ver2.1}) (applied to $S\setminus\left\{  t\right\}  $
instead of $S$) shows that $S\setminus\left\{  t\right\}  $ has a maximum.
Also, Lemma \ref{lem.ind.max-1el} (applied to $x=t$) shows that the set
$\left\{  t\right\}  $ has a maximum, namely $t$. Hence, Proposition
\ref{prop.ind.max-PuQ} (applied to $P=S\setminus\left\{  t\right\}  $ and
$Q=\left\{  t\right\}  $) shows that the set $\left(  S\setminus\left\{
t\right\}  \right)  \cup\left\{  t\right\}  $ has a maximum. Since $\left(
S\setminus\left\{  t\right\}  \right)  \cup\left\{  t\right\}  =S$, this
rewrites as follows: The set $S$ has a maximum.

Now, forget that we fixed $S$. We thus have shown that if $S$ is a nonempty
finite set of integers satisfying $\left\vert S\right\vert -1=m+1$, then $S$
has a maximum. In other words, Theorem \ref{thm.ind.max} holds under the
condition that $\left\vert S\right\vert -1=m+1$. This completes the induction
step. Thus, the induction proof of Theorem \ref{thm.ind.max} is complete.
\end{proof}

We could have shortened this proof even further if we didn't explicitly state
(\ref{pf.thm.ind.max.ver2.1}), but rather (instead of applying
(\ref{pf.thm.ind.max.ver2.1})) said that \textquotedblleft we can apply
Theorem \ref{thm.ind.max} to $S\setminus\left\{  t\right\}  $ instead of
$S$\textquotedblright.

Let us stress again that, in order to prove Theorem \ref{thm.ind.max} by
induction on $\left\vert S\right\vert -1$, we had to check that $\left\vert
S\right\vert -1\in\mathbb{N}$ whenever $S$ satisfies the assumptions of
Theorem \ref{thm.ind.max}.\footnote{In our first version of the above proof,
we checked this at the end; in the second version, we checked it at the
beginning of the proof.} This check was necessary. For example, if we had
instead tried to proceed by induction on $\left\vert S\right\vert -2$, then we
would only have proven Theorem \ref{thm.ind.max} under the condition that
$\left\vert S\right\vert -2\in\mathbb{N}$; but this condition isn't always
satisfied (indeed, it misses the case when $S$ is a $1$-element set).

\subsubsection{Vacuous truth and induction bases}

Can we also prove Theorem \ref{thm.ind.max} by induction on $\left\vert
S\right\vert $ (instead of $\left\vert S\right\vert -1$)? This seems a bit
strange, since $\left\vert S\right\vert $ can never be $0$ in Theorem
\ref{thm.ind.max} (because $S$ is required to be nonempty), so that the
induction base would be talking about a situation that never occurs. However,
there is nothing wrong about it, and we already do talk about such situations
oftentimes (for example, every time we make a proof by contradiction). The
following concept from basic logic explains this:

\begin{convention}
\label{conv.logic.vacuous} \textbf{(a)} A logical statement of the form
\textquotedblleft if $\mathcal{A}$, then $\mathcal{B}$\textquotedblright%
\ (where $\mathcal{A}$ and $\mathcal{B}$ are two statements) is said to be
\textit{vacuously true} if $\mathcal{A}$ does not hold. For example, the
statement \textquotedblleft if $0=1$, then every set is
empty\textquotedblright\ is vacuously true, because $0=1$ is false. The
statement \textquotedblleft if $0=1$, then $1=1$\textquotedblright\ is also
vacuously true, although its truth can also be seen as a consequence of the
fact that $1=1$ is true.

By the laws of logic, a vacuously true statement is always true! This may
sound counterintuitive, but actually makes perfect sense: A statement
\textquotedblleft if $\mathcal{A}$, then $\mathcal{B}$\textquotedblright\ only
says anything about situations where $\mathcal{A}$ holds. If $\mathcal{A}$
never holds, then it therefore says nothing. And when you are saying nothing,
you are certainly not lying.

The principle that a vacuously true statement always holds is known as
\textquotedblleft\textit{ex falso quodlibet}\textquotedblright\ (literal
translation: \textquotedblleft from the false, anything\textquotedblright) or
\textquotedblleft%
\href{https://en.wikipedia.org/wiki/Principle_of_explosion}{\textit{principle
of explosion}}\textquotedblright. It can be restated as follows: From a false
statement, any statement follows.

\textbf{(b)} Now, let $X$ be a set, and let $\mathcal{A}\left(  x\right)  $
and $\mathcal{B}\left(  x\right)  $ be two statements defined for each $x\in
X$. A statement of the form \textquotedblleft for each $x\in X$ satisfying
$\mathcal{A}\left(  x\right)  $, we have $\mathcal{B}\left(  x\right)
$\textquotedblright\ will automatically hold if there exists no $x\in X$
satisfying $\mathcal{A}\left(  x\right)  $. (Indeed, this statement can be
rewritten as \textquotedblleft for each $x\in X$, we have $\left(  \text{if
}\mathcal{A}\left(  x\right)  \text{, then }\mathcal{B}\left(  x\right)
\right)  $\textquotedblright; but this holds because the statement
\textquotedblleft if $\mathcal{A}\left(  x\right)  $, then $\mathcal{B}\left(
x\right)  $\textquotedblright\ is vacuously true for each $x\in X$.) Such a
statement will also be called \textit{vacuously true}.

For example, the statement \textquotedblleft if $n\in\mathbb{N}$ is both odd
and even, then $n=n+1$\textquotedblright\ is vacuously true, since no
$n\in\mathbb{N}$ can be both odd and even at the same time.

\textbf{(c)} Now, let $X$ be the empty set (that is, $X=\varnothing$), and let
$\mathcal{B}\left(  x\right)  $ be a statement defined for each $x\in X$.
Then, a statement of the form \textquotedblleft for each $x\in X$, we have
$\mathcal{B}\left(  x\right)  $\textquotedblright\ will automatically hold.
(Indeed, this statement can be rewritten as \textquotedblleft for each $x\in
X$, we have $\left(  \text{if }x\in X\text{, then }\mathcal{B}\left(
x\right)  \right)  $\textquotedblright; but this holds because the statement
\textquotedblleft if $x\in X$, then $\mathcal{B}\left(  x\right)
$\textquotedblright\ is vacuously true for each $x\in X$, since its premise
($x\in X$) is false.) Again, such a statement is said to be \textit{vacuously
true}.

For example, the statement \textquotedblleft for each $x\in\varnothing$, we
have $x\neq x$\textquotedblright\ is vacuously true (because there exists no
$x\in\varnothing$).
\end{convention}

Thus, if we try to prove Theorem \ref{thm.ind.max} by induction on $\left\vert
S\right\vert $, then the induction base becomes vacuously true. However, the
induction step becomes more complicated, since we can no longer argue that
$S\setminus\left\{  t\right\}  $ is nonempty, but instead have to account for
the case when $S\setminus\left\{  t\right\}  $ is empty as well. So we gain
and we lose at the same time. Here is how this proof looks like:

\begin{proof}
[Second proof of Theorem \ref{thm.ind.max}.]Clearly, $\left\vert S\right\vert
\in\mathbb{N}$ (since $S$ is a finite set). Hence, we can apply induction on
$\left\vert S\right\vert $ to prove Theorem \ref{thm.ind.max}:

\textit{Induction base:} Theorem \ref{thm.ind.max} holds under the condition
that $\left\vert S\right\vert =0$\ \ \ \ \footnote{\textit{Proof.} Let $S$ be
as in Theorem \ref{thm.ind.max}, and assume that $\left\vert S\right\vert =0$.
We must show that the claim of Theorem \ref{thm.ind.max} holds.
\par
Indeed, $\left\vert S\right\vert =0$, so that $S$ is the empty set. This
contradicts the assumption that $S$ be nonempty. From this contradiction, we
conclude that everything holds (by the \textquotedblleft ex falso
quodlibet\textquotedblright\ principle). Thus, in particular, the claim of
Theorem \ref{thm.ind.max} holds. This completes our proof.}. This completes
the induction base.

\textit{Induction step:} Let $m\in\mathbb{N}$. Assume that Theorem
\ref{thm.ind.max} holds under the condition that $\left\vert S\right\vert =m$.
We shall now show that Theorem \ref{thm.ind.max} holds under the condition
that $\left\vert S\right\vert =m+1$.

We have assumed that Theorem \ref{thm.ind.max} holds under the condition that
$\left\vert S\right\vert =m$. In other words,%
\begin{equation}
\left(
\begin{array}
[c]{c}%
\text{if }S\text{ is a nonempty finite set of integers satisfying }\left\vert
S\right\vert =m\text{,}\\
\text{then }S\text{ has a maximum}%
\end{array}
\right)  . \label{pf.thm.ind.max.ver3.1}%
\end{equation}

Now, let $S$ be a nonempty finite set of integers satisfying $\left\vert
S\right\vert =m+1$. We want to prove that $S$ has a maximum.

There exists some $t\in S$ (since $S$ is nonempty). Consider this $t$. We have
$\left(  S\setminus\left\{  t\right\}  \right)  \cup\left\{  t\right\}
=S\cup\left\{  t\right\}  =S$ (since $t\in S$). Lemma \ref{lem.ind.max-1el}
(applied to $x=t$) shows that the set $\left\{  t\right\}  $ has a maximum,
namely $t$.

We are in one of the following two cases:

\textit{Case 1:} We have $S\setminus\left\{  t\right\}  =\varnothing$.

\textit{Case 2:} We have $S\setminus\left\{  t\right\}  \neq\varnothing$.

Let us first consider Case 1. In this case, we have $S\setminus\left\{
t\right\}  =\varnothing$. Hence, $S\subseteq\left\{  t\right\}  $. Thus,
either $S=\varnothing$ or $S=\left\{  t\right\}  $ (since the only subsets of
$\left\{  t\right\}  $ are $\varnothing$ and $\left\{  t\right\}  $). Since
$S=\varnothing$ is impossible (because $S$ is nonempty), we thus have
$S=\left\{  t\right\}  $. But the set $\left\{  t\right\}  $ has a maximum. In
view of $S=\left\{  t\right\}  $, this rewrites as follows: The set $S$ has a
maximum. Thus, our goal (to prove that $S$ has a maximum) is achieved in Case 1.

Let us now consider Case 2. In this case, we have $S\setminus\left\{
t\right\}  \neq\varnothing$. Hence, the set $S\setminus\left\{  t\right\}  $
is nonempty. From $t\in S$, we obtain $\left\vert S\setminus\left\{
t\right\}  \right\vert =\left\vert S\right\vert -1=m$ (since $\left\vert
S\right\vert =m+1$). Furthermore, the set $S\setminus\left\{  t\right\}  $ is
finite (since $S$ is finite). Hence, (\ref{pf.thm.ind.max.ver3.1}) (applied to
$S\setminus\left\{  t\right\}  $ instead of $S$) shows that $S\setminus
\left\{  t\right\}  $ has a maximum. Also, recall that the set $\left\{
t\right\}  $ has a maximum. Hence, Proposition \ref{prop.ind.max-PuQ} (applied
to $P=S\setminus\left\{  t\right\}  $ and $Q=\left\{  t\right\}  $) shows that
the set $\left(  S\setminus\left\{  t\right\}  \right)  \cup\left\{
t\right\}  $ has a maximum. Since $\left(  S\setminus\left\{  t\right\}
\right)  \cup\left\{  t\right\}  =S$, this rewrites as follows: The set $S$
has a maximum. Hence, our goal (to prove that $S$ has a maximum) is achieved
in Case 2.

We have now proven that $S$ has a maximum in each of the two Cases 1 and 2.
Therefore, $S$ always has a maximum (since Cases 1 and 2 cover all possibilities).

Now, forget that we fixed $S$. We thus have shown that if $S$ is a nonempty
finite set of integers satisfying $\left\vert S\right\vert =m+1$, then $S$ has
a maximum. In other words, Theorem \ref{thm.ind.max} holds under the condition
that $\left\vert S\right\vert =m+1$. This completes the induction step. Thus,
the induction proof of Theorem \ref{thm.ind.max} is complete.
\end{proof}

\subsubsection{Further results on maxima and minima}

We can replace \textquotedblleft integers\textquotedblright\ by
\textquotedblleft rational numbers\textquotedblright\ or \textquotedblleft
real numbers\textquotedblright\ in Theorem \ref{thm.ind.max}; all the proofs
given above still apply then. Thus, we obtain the following:

\begin{theorem}
\label{thm.ind.max2}Let $S$ be a nonempty finite set of integers (or rational
numbers, or real numbers). Then, $S$ has a maximum.
\end{theorem}

Hence, if $S$ is a nonempty finite set of integers (or rational numbers, or
real numbers), then $\max S$ is well-defined (because Theorem
\ref{thm.ind.max2} shows that $S$ has a maximum, and Proposition
\ref{prop.ind.max-uni} shows that this maximum is unique).

Moreover, just as we have defined maxima (i.e., largest elements) of sets, we
can define minima (i.e., smallest elements) of sets, and prove similar results
about them:

\begin{definition}
Let $S$ be a set of integers (or rational numbers, or real numbers). A
\textit{minimum} of $S$ is defined to be an element $s\in S$ that satisfies%
\[
\left(  s\leq t\text{ for each }t\in S\right)  .
\]
In other words, a minimum of $S$ is defined to be an element of $S$ which is
less or equal to each element of $S$.

(The plural of the word \textquotedblleft minimum\textquotedblright\ is
\textquotedblleft minima\textquotedblright.)
\end{definition}

\begin{example}
The set $\left\{  2,4,5\right\}  $ has exactly one minimum: namely, $2$.

The set $\mathbb{N}=\left\{  0,1,2,\ldots\right\}  $ has exactly one minimum:
namely, $0$.

The set $\left\{  0,-1,-2,\ldots\right\}  $ has no minimum: If $k$ was a
minimum of this set, then we would have $k\leq k-1$, which is absurd.

The set $\varnothing$ has no minimum, since a minimum would have to be an
element of $\varnothing$.
\end{example}

The analogue of Proposition \ref{prop.ind.max-uni} for minima instead of
maxima looks exactly as one would expect it:

\begin{proposition}
\label{prop.ind.min-uni}Let $S$ be a set of integers (or rational numbers, or
real numbers). Then, $S$ has \textbf{at most one} minimum.
\end{proposition}

\begin{proof}
[Proof of Proposition \ref{prop.ind.min-uni}.]To obtain a proof of Proposition
\ref{prop.ind.min-uni}, it suffices to replace every \textquotedblleft$\geq
$\textquotedblright\ sign by a \textquotedblleft$\leq$\textquotedblright\ sign
(and every word \textquotedblleft maximum\textquotedblright\ by
\textquotedblleft minimum\textquotedblright) in the proof of Proposition
\ref{prop.ind.max-uni} given above.
\end{proof}

\begin{definition}
Let $S$ be a set of integers (or rational numbers, or real numbers).
Proposition \ref{prop.ind.min-uni} shows that $S$ has \textbf{at most one}
minimum. Thus, if $S$ has a minimum, then this minimum is the unique minimum
of $S$; we shall thus call it \textit{the minimum} of $S$ or \textit{the
smallest element} of $S$. We shall denote this minimum by $\min S$.
\end{definition}

The analogue of Theorem \ref{thm.ind.max2} is the following:

\begin{theorem}
\label{thm.ind.min2}Let $S$ be a nonempty finite set of integers (or rational
numbers, or real numbers). Then, $S$ has a minimum.
\end{theorem}

\begin{proof}
[Proof of Theorem \ref{thm.ind.min2}.]To obtain a proof of Theorem
\ref{thm.ind.min2}, it suffices to replace every \textquotedblleft$\geq
$\textquotedblright\ sign by a \textquotedblleft$\leq$\textquotedblright\ sign
(and every word \textquotedblleft maximum\textquotedblright\ by
\textquotedblleft minimum\textquotedblright) in the proof of Theorem
\ref{thm.ind.max2} given above (and also in the proofs of all the auxiliary
results that were used in said proof).\footnote{To be technically precise: not
every \textquotedblleft$\geq$\textquotedblright\ sign, of course. The
\textquotedblleft$\geq$\textquotedblright\ sign in \textquotedblleft$m\geq
0$\textquotedblright\ should stay unchanged.}
\end{proof}

Alternatively, Theorem \ref{thm.ind.min2} can be obtained from Theorem
\ref{thm.ind.max2} by applying the latter theorem to the set $\left\{
-s\ \mid\ s\in S\right\}  $. In fact, it is easy to see that a number $x$ is
the minimum of $S$ if and only if $-x$ is the maximum of the set $\left\{
-s\ \mid\ s\in S\right\}  $. We leave the details of this simple argument to
the reader.

We also should mention that Theorem \ref{thm.ind.min2} holds \textbf{without}
requiring that $S$ be finite, if we instead require that $S$ consist of
nonnegative integers:

\begin{theorem}
\label{thm.ind.max3}Let $S$ be a nonempty set of nonnegative integers. Then,
$S$ has a minimum.
\end{theorem}

But $S$ does not necessarily have a maximum in this situation; the
nonnegativity requirement has \textquotedblleft broken the
symmetry\textquotedblright\ between maxima and minima.

We note that the word \textquotedblleft integers\textquotedblright\ is crucial
in Theorem \ref{thm.ind.max3}. If we replaced \textquotedblleft
integers\textquotedblright\ by \textquotedblleft rational
numbers\textquotedblright, then the theorem would no longer hold (for example,
the set of all positive rational numbers has no minimum, since positive
rational numbers can get arbitrarily close to $0$ yet cannot equal $0$).

\begin{proof}
[Proof of Theorem \ref{thm.ind.max3}.]The set $S$ is nonempty. Thus, there
exists some $p\in S$. Consider this $p$.

We have $p\in S\subseteq\mathbb{N}$ (since $S$ is a set of nonnegative
integers). Thus, $p\in\left\{  0,1,\ldots,p\right\}  $. Combining this with
$p\in S$, we obtain $p\in\left\{  0,1,\ldots,p\right\}  \cap S$. Hence, the
set $\left\{  0,1,\ldots,p\right\}  \cap S$ contains the element $p$, and thus
is nonempty. Moreover, this set $\left\{  0,1,\ldots,p\right\}  \cap S$ is a
subset of the finite set $\left\{  0,1,\ldots,p\right\}  $, and thus is finite.

Now we know that $\left\{  0,1,\ldots,p\right\}  \cap S$ is a nonempty finite
set of integers. Hence, Theorem \ref{thm.ind.min2} (applied to $\left\{
0,1,\ldots,p\right\}  \cap S$ instead of $S$) shows that the set $\left\{
0,1,\ldots,p\right\}  \cap S$ has a minimum. Denote this minimum by $m$.

Hence, $m$ is a minimum of the set $\left\{  0,1,\ldots,p\right\}  \cap S$. In
other words, $m$ is an element $s\in\left\{  0,1,\ldots,p\right\}  \cap S$
that satisfies%
\[
\left(  s\leq t\text{ for each }t\in\left\{  0,1,\ldots,p\right\}  \cap
S\right)
\]
(by the definition of a minimum). In other words, $m$ is an element of
$\left\{  0,1,\ldots,p\right\}  \cap S$ and satisfies%
\begin{equation}
\left(  m\leq t\text{ for each }t\in\left\{  0,1,\ldots,p\right\}  \cap
S\right)  . \label{pf.thm.ind.max3.m1}%
\end{equation}

Hence, $m\in\left\{  0,1,\ldots,p\right\}  \cap S\subseteq\left\{
0,1,\ldots,p\right\}  $, so that $m\leq p$.

Furthermore, $m\in\left\{  0,1,\ldots,p\right\}  \cap S\subseteq S$. Moreover,
we have%
\begin{equation}
\left(  m\leq t\text{ for each }t\in S\right)  . \label{pf.thm.ind.max3.m2}%
\end{equation}

\begin{vershort}
[\textit{Proof of (\ref{pf.thm.ind.max3.m2}):} Let $t\in S$. We must prove
that $m\leq t$.

If $t\in\left\{  0,1,\ldots,p\right\}  \cap S$, then this follows from
(\ref{pf.thm.ind.max3.m1}). Hence, for the rest of this proof, we can WLOG
assume that we don't have $t\in\left\{  0,1,\ldots,p\right\}  \cap S$. Assume
this. Thus, $t\notin\left\{  0,1,\ldots,p\right\}  \cap S$. Combining $t\in S$
with $t\notin\left\{  0,1,\ldots,p\right\}  \cap S$, we obtain
\[
t\in S\setminus\left(  \left\{  0,1,\ldots,p\right\}  \cap S\right)
=S\setminus\left\{  0,1,\ldots,p\right\}  .
\]
Hence, $t\notin\left\{  0,1,\ldots,p\right\}  $, so that $t>p$ (since
$t\in\mathbb{N}$). Therefore, $t\geq p\geq m$ (since $m\leq p$), so that
$m\leq t$. This completes the proof of (\ref{pf.thm.ind.max3.m2}).]
\end{vershort}

\begin{verlong}
[\textit{Proof of (\ref{pf.thm.ind.max3.m2}):} Let $t\in S$. We must prove
that $m\leq t$.

If $t\in\left\{  0,1,\ldots,p\right\}  \cap S$, then this follows from
(\ref{pf.thm.ind.max3.m1}). Hence, for the rest of this proof, we can WLOG
assume that we don't have $t\in\left\{  0,1,\ldots,p\right\}  \cap S$. Assume
this. Thus, $t\notin\left\{  0,1,\ldots,p\right\}  \cap S$ (since we don't
have $t\in\left\{  0,1,\ldots,p\right\}  \cap S$). Combining $t\in S$ with
$t\notin\left\{  0,1,\ldots,p\right\}  \cap S$, we obtain
\begin{align*}
t  &  \in S\setminus\left(  \left\{  0,1,\ldots,p\right\}  \cap S\right)
=\underbrace{S}_{\subseteq\mathbb{N}}\setminus\left\{  0,1,\ldots,p\right\} \\
&  \subseteq\mathbb{N}\setminus\left\{  0,1,\ldots,p\right\}  =\left\{
p+1,p+2,p+3,\ldots\right\}  .
\end{align*}
Hence, $t\geq p+1\geq p\geq m$ (since $m\leq p$), so that $m\leq t$. This
completes the proof of (\ref{pf.thm.ind.max3.m2}).]
\end{verlong}

Now, we know that $m$ is an element of $S$ (since $m\in S$) and satisfies
\newline$\left(  m\leq t\text{ for each }t\in S\right)  $ (by
(\ref{pf.thm.ind.max3.m2})). In other words, $m$ is an $s\in S$ that satisfies
\newline$\left(  s\leq t\text{ for each }t\in S\right)  $. In other words, $m$
is a minimum of $S$ (by the definition of a minimum). Thus, $S$ has a minimum
(namely, $m$). This proves Theorem \ref{thm.ind.max3}.
\end{proof}

\subsection{Increasing lists of finite sets}

We shall next study (again using induction) another basic feature of finite sets.

We recall that \textquotedblleft list\textquotedblright\ is just a synonym for
\textquotedblleft tuple\textquotedblright; i.e., a list is a $k$-tuple for
some $k\in\mathbb{N}$. Note that tuples and lists are always understood to be
finite and ordered.

\begin{definition}
\label{def.ind.inclist0}Let $S$ be a set of integers. An \textit{increasing
list} of $S$ shall mean a list $\left(  s_{1},s_{2},\ldots,s_{k}\right)  $ of
elements of $S$ such that $S=\left\{  s_{1},s_{2},\ldots,s_{k}\right\}  $ and
$s_{1}<s_{2}<\cdots<s_{k}$.
\end{definition}

In other words, if $S$ is a set of integers, then an increasing list of $S$
means a list such that

\begin{itemize}
\item the set $S$ consists of all elements of this list, and

\item the elements of this list are strictly increasing.
\end{itemize}

For example, $\left(  2,4,6\right)  $ is an increasing list of the set
$\left\{  2,4,6\right\}  $, but neither $\left(  2,6\right)  $ nor $\left(
2,4,4,6\right)  $ nor $\left(  4,2,6\right)  $ nor $\left(  2,4,5,6\right)  $
is an increasing list of this set. For another example, $\left(
1,4,9,16\right)  $ is an increasing list of the set $\left\{  i^{2}%
\ \mid\ i\in\left\{  1,2,3,4\right\}  \right\}  =\left\{  1,4,9,16\right\}  $.
For yet another example, the empty list $\left(  {}\right)  $ is an increasing
list of the empty set $\varnothing$.

Now, it is intuitively obvious that any finite set $S$ of integers has a
unique increasing list -- we just need to list all the elements of $S$ in
increasing order, with no repetitions. But from the viewpoint of rigorous
mathematics, this needs to be proven. Let us state this as a theorem:

\begin{theorem}
\label{thm.ind.inclist.unex}Let $S$ be a finite set of integers. Then, $S$ has
exactly one increasing list.
\end{theorem}

Before we prove this theorem, let us show some auxiliary facts:

\begin{proposition}
\label{prop.ind.inclist.size}Let $S$ be a set of integers. Let $\left(
s_{1},s_{2},\ldots,s_{k}\right)  $ be an increasing list of $S$. Then:

\textbf{(a)} The set $S$ is finite.

\textbf{(b)} We have $\left\vert S\right\vert =k$.

\textbf{(c)} The elements $s_{1},s_{2},\ldots,s_{k}$ are distinct.
\end{proposition}

\begin{proof}
[Proof of Proposition \ref{prop.ind.inclist.size}.]We know that $\left(
s_{1},s_{2},\ldots,s_{k}\right)  $ is an increasing list of $S$. In other
words, $\left(  s_{1},s_{2},\ldots,s_{k}\right)  $ is a list of elements of
$S$ such that $S=\left\{  s_{1},s_{2},\ldots,s_{k}\right\}  $ and $s_{1}%
<s_{2}<\cdots<s_{k}$ (by the definition of an \textquotedblleft increasing
list\textquotedblright).

From $S=\left\{  s_{1},s_{2},\ldots,s_{k}\right\}  $, we conclude that the set
$S$ has at most $k$ elements. Thus, the set $S$ is finite. This proves
Proposition \ref{prop.ind.inclist.size} \textbf{(a)}.

We have $s_{1}<s_{2}<\cdots<s_{k}$. Hence, if $u$ and $v$ are two elements of
$\left\{  1,2,\ldots,k\right\}  $ such that $u<v$, then $s_{u}<s_{v}$ (by
Corollary \ref{cor.mod.chain-ineq3}, applied to $a_{i}=s_{i}$) and therefore
$s_{u}\neq s_{v}$. In other words, the elements $s_{1},s_{2},\ldots,s_{k}$ are
distinct. This proves Proposition \ref{prop.ind.inclist.size} \textbf{(c)}.

The $k$ elements $s_{1},s_{2},\ldots,s_{k}$ are distinct; thus, the set
$\left\{  s_{1},s_{2},\ldots,s_{k}\right\}  $ has size $k$. In other words,
the set $S$ has size $k$ (since $S=\left\{  s_{1},s_{2},\ldots,s_{k}\right\}
$). In other words, $\left\vert S\right\vert =k$. This proves Proposition
\ref{prop.ind.inclist.size} \textbf{(b)}.
\end{proof}

\begin{proposition}
\label{prop.ind.inclist.empty}The set $\varnothing$ has exactly one increasing
list: namely, the empty list $\left(  {}\right)  $.
\end{proposition}

\begin{proof}
[Proof of Proposition \ref{prop.ind.inclist.empty}.]The empty list $\left(
{}\right)  $ satisfies $\varnothing=\left\{  {}\right\}  $. Thus, the empty
list $\left(  {}\right)  $ is a list $\left(  s_{1},s_{2},\ldots,s_{k}\right)
$ of elements of $\varnothing$ such that $\varnothing=\left\{  s_{1}%
,s_{2},\ldots,s_{k}\right\}  $ and $s_{1}<s_{2}<\cdots<s_{k}$ (indeed, the
chain of inequalities $s_{1}<s_{2}<\cdots<s_{k}$ is vacuously true for the
empty list $\left(  {}\right)  $, because it contains no inequality signs). In
other words, the empty list $\left(  {}\right)  $ is an increasing list of
$\varnothing$ (by the definition of an increasing list). It remains to show
that it is the only increasing list of $\varnothing$.

Let $\left(  s_{1},s_{2},\ldots,s_{k}\right)  $ be any increasing list of
$\varnothing$. Then, Proposition \ref{prop.ind.inclist.size} \textbf{(b)}
(applied to $S=\varnothing$) yields $\left\vert \varnothing\right\vert =k$.
Hence, $k=\left\vert \varnothing\right\vert =0$, so that $\left(  s_{1}%
,s_{2},\ldots,s_{k}\right)  =\left(  s_{1},s_{2},\ldots,s_{0}\right)  =\left(
{}\right)  $.

Now, forget that we fixed $\left(  s_{1},s_{2},\ldots,s_{k}\right)  $. We thus
have shown that if $\left(  s_{1},s_{2},\ldots,s_{k}\right)  $ is any
increasing list of $\varnothing$, then $\left(  s_{1},s_{2},\ldots
,s_{k}\right)  =\left(  {}\right)  $. In other words, any increasing list of
$\varnothing$ is $\left(  {}\right)  $. Therefore, the set $\varnothing$ has
exactly one increasing list: namely, the empty list $\left(  {}\right)  $
(since we already know that $\left(  {}\right)  $ is an increasing list of
$\varnothing$). This proves Proposition \ref{prop.ind.inclist.empty}.
\end{proof}

\begin{proposition}
\label{prop.ind.inclist.nonempty1}Let $S$ be a nonempty finite set of
integers. Let $m=\max S$. Let $\left(  s_{1},s_{2},\ldots,s_{k}\right)  $ be
any increasing list of $S$. Then:

\textbf{(a)} We have $k\geq1$ and $s_{k}=m$.

\textbf{(b)} The list $\left(  s_{1},s_{2},\ldots,s_{k-1}\right)  $ is an
increasing list of $S\setminus\left\{  m\right\}  $.
\end{proposition}

\begin{proof}
[Proof of Proposition \ref{prop.ind.inclist.nonempty1}.]We know that $\left(
s_{1},s_{2},\ldots,s_{k}\right)  $ is an increasing list of $S$. In other
words, $\left(  s_{1},s_{2},\ldots,s_{k}\right)  $ is a list of elements of
$S$ such that $S=\left\{  s_{1},s_{2},\ldots,s_{k}\right\}  $ and $s_{1}%
<s_{2}<\cdots<s_{k}$ (by the definition of an \textquotedblleft increasing
list\textquotedblright).

Proposition \ref{prop.ind.inclist.size} \textbf{(b)} yields $\left\vert
S\right\vert =k$. Hence, $k=\left\vert S\right\vert >0$ (since $S$ is
nonempty). Thus, $k\geq1$ (since $k$ is an integer). Therefore, $s_{k}$ is
well-defined. Clearly, $k\in\left\{  1,2,\ldots,k\right\}  $ (since $k\geq1$),
so that $s_{k}\in\left\{  s_{1},s_{2},\ldots,s_{k}\right\}  =S$.

We have $s_{1}<s_{2}<\cdots<s_{k}$ and thus $s_{1}\leq s_{2}\leq\cdots\leq
s_{k}$. Hence, Proposition \ref{prop.mod.chain-ineq} (applied to $a_{i}=s_{i}%
$) shows that if $u$ and $v$ are two elements of $\left\{  1,2,\ldots
,k\right\}  $ such that $u\leq v$, then%
\begin{equation}
s_{u}\leq s_{v}. \label{pf.prop.ind.inclist.nonempty1.1}%
\end{equation}

Thus, we have $\left(  s_{k}\geq t\text{ for each }t\in S\right)
$\ \ \ \ \footnote{\textit{Proof.} Let $t\in S$. Thus, $t\in S=\left\{
s_{1},s_{2},\ldots,s_{k}\right\}  $. Hence, $t=s_{u}$ for some $u\in\left\{
1,2,\ldots,k\right\}  $. Consider this $u$. Now, $u$ and $k$ are elements of
$\left\{  1,2,\ldots,k\right\}  $ such that $u\leq k$ (since $u\in\left\{
1,2,\ldots,k\right\}  $). Hence, (\ref{pf.prop.ind.inclist.nonempty1.1})
(applied to $v=k$) yields $s_{u}\leq s_{k}$. Hence, $s_{k}\geq s_{u}=t$ (since
$t=s_{u}$), qed.}. Hence, $s_{k}$ is an element $s\in S$ that satisfies
$\left(  s\geq t\text{ for each }t\in S\right)  $ (since $s_{k}\in S$). In
other words, $s_{k}$ is a maximum of $S$ (by the definition of a maximum).
Since we know that $S$ has at most one maximum (by Proposition
\ref{prop.ind.max-uni}), we thus conclude that $s_{k}$ is \textbf{the} maximum
of $S$. In other words, $s_{k}=\max S$. Hence, $s_{k}=\max S=m$. This
completes the proof of Proposition \ref{prop.ind.inclist.nonempty1}
\textbf{(a)}.

\textbf{(b)} From $s_{1}<s_{2}<\cdots<s_{k}$, we obtain $s_{1}<s_{2}%
<\cdots<s_{k-1}$. Furthermore, the elements $s_{1},s_{2},\ldots,s_{k}$ are
distinct (according to Proposition \ref{prop.ind.inclist.size} \textbf{(c)}).
In other words, for any two distinct elements $u$ and $v$ of $\left\{
1,2,\ldots,k\right\}  $, we have%
\begin{equation}
s_{u}\neq s_{v}. \label{pf.prop.ind.inclist.nonempty1.2}%
\end{equation}
Hence, $s_{k}\notin\left\{  s_{1},s_{2},\ldots,s_{k-1}\right\}  $%
\ \ \ \ \footnote{\textit{Proof.} Assume the contrary. Thus, $s_{k}\in\left\{
s_{1},s_{2},\ldots,s_{k-1}\right\}  $. In other words, $s_{k}=s_{u}$ for some
$u\in\left\{  1,2,\ldots,k-1\right\}  $. Consider this $u$. We have
$u\in\left\{  1,2,\ldots,k-1\right\}  \subseteq\left\{  1,2,\ldots,k\right\}
$.
\par
Now, $u\in\left\{  1,2,\ldots,k-1\right\}  $, so that $u\leq k-1<k$ and thus
$u\neq k$. Hence, the elements $u$ and $k$ of $\left\{  1,2,\ldots,k\right\}
$ are distinct. Thus, (\ref{pf.prop.ind.inclist.nonempty1.2}) (applied to
$v=k$) yields $s_{u}\neq s_{k}=s_{u}$. This is absurd. This contradiction
shows that our assumption was wrong, qed.}. Now,%
\begin{align*}
\underbrace{S}_{\substack{=\left\{  s_{1},s_{2},\ldots,s_{k}\right\}
\\=\left\{  s_{1},s_{2},\ldots,s_{k-1}\right\}  \cup\left\{  s_{k}\right\}
}}\setminus\left\{  \underbrace{m}_{=s_{k}}\right\}   &  =\left(  \left\{
s_{1},s_{2},\ldots,s_{k-1}\right\}  \cup\left\{  s_{k}\right\}  \right)
\setminus\left\{  s_{k}\right\} \\
&  =\left\{  s_{1},s_{2},\ldots,s_{k-1}\right\}  \setminus\left\{
s_{k}\right\}  =\left\{  s_{1},s_{2},\ldots,s_{k-1}\right\}
\end{align*}
(since $s_{k}\notin\left\{  s_{1},s_{2},\ldots,s_{k-1}\right\}  $). Hence, the
elements $s_{1},s_{2},\ldots,s_{k-1}$ belong to the set $S\setminus\left\{
m\right\}  $ (since they clearly belong to the set $\left\{  s_{1}%
,s_{2},\ldots,s_{k-1}\right\}  =S\setminus\left\{  m\right\}  $). In other
words, $\left(  s_{1},s_{2},\ldots,s_{k-1}\right)  $ is a list of elements of
$S\setminus\left\{  m\right\}  $.

Now, we know that $\left(  s_{1},s_{2},\ldots,s_{k-1}\right)  $ is a list of
elements of $S\setminus\left\{  m\right\}  $ such that $S\setminus\left\{
m\right\}  =\left\{  s_{1},s_{2},\ldots,s_{k-1}\right\}  $ and $s_{1}%
<s_{2}<\cdots<s_{k-1}$. In other words, $\left(  s_{1},s_{2},\ldots
,s_{k-1}\right)  $ is an increasing list of $S\setminus\left\{  m\right\}  $.
This proves Proposition \ref{prop.ind.inclist.nonempty1} \textbf{(b)}.
\end{proof}

We are now ready to prove Theorem \ref{thm.ind.inclist.unex}:

\begin{proof}
[Proof of Theorem \ref{thm.ind.inclist.unex}.]We shall prove Theorem
\ref{thm.ind.inclist.unex} by induction on $\left\vert S\right\vert $:

\textit{Induction base:} Theorem \ref{thm.ind.inclist.unex} holds under the
condition that $\left\vert S\right\vert =0$\ \ \ \ \footnote{\textit{Proof.}
Let $S$ be as in Theorem \ref{thm.ind.inclist.unex}, and assume that
$\left\vert S\right\vert =0$. We must show that the claim of Theorem
\ref{thm.ind.inclist.unex} holds.
\par
Indeed, $\left\vert S\right\vert =0$, so that $S$ is the empty set. Thus,
$S=\varnothing$. But Proposition \ref{prop.ind.inclist.empty} shows that the
set $\varnothing$ has exactly one increasing list. Since $S=\varnothing$, this
rewrites as follows: The set $S$ has exactly one increasing list. Thus, the
claim of Theorem \ref{thm.ind.inclist.unex} holds. This completes our proof.}.
This completes the induction base.

\textit{Induction step:} Let $g\in\mathbb{N}$. Assume that Theorem
\ref{thm.ind.inclist.unex} holds under the condition that $\left\vert
S\right\vert =g$. We shall now show that Theorem \ref{thm.ind.inclist.unex}
holds under the condition that $\left\vert S\right\vert =g+1$.

We have assumed that Theorem \ref{thm.ind.inclist.unex} holds under the
condition that $\left\vert S\right\vert =g$. In other words,%
\begin{equation}
\left(
\begin{array}
[c]{c}%
\text{if }S\text{ is a finite set of integers satisfying }\left\vert
S\right\vert =g\text{,}\\
\text{then }S\text{ has exactly one increasing list}%
\end{array}
\right)  . \label{pf.thm.ind.inclist.unex.IH}%
\end{equation}

Now, let $S$ be a finite set of integers satisfying $\left\vert S\right\vert
=g+1$. We want to prove that $S$ has exactly one increasing list.

The set $S$ is nonempty (since $\left\vert S\right\vert =g+1>g\geq0$). Thus,
$S$ has a maximum (by Theorem \ref{thm.ind.max}). Hence, $\max S$ is
well-defined. Set $m=\max S$. Thus, $m=\max S\in S$ (by
(\ref{eq.ind.max.def-max.1})). Therefore, $\left\vert S\setminus\left\{
m\right\}  \right\vert =\left\vert S\right\vert -1=g$ (since $\left\vert
S\right\vert =g+1$). Hence, (\ref{pf.thm.ind.inclist.unex.IH}) (applied to
$S\setminus\left\{  m\right\}  $ instead of $S$) shows that $S\setminus
\left\{  m\right\}  $ has exactly one increasing list. Let $\left(
t_{1},t_{2},\ldots,t_{j}\right)  $ be this list. We extend this list to a
$\left(  j+1\right)  $-tuple $\left(  t_{1},t_{2},\ldots,t_{j+1}\right)  $ by
setting $t_{j+1}=m$.

We have defined $\left(  t_{1},t_{2},\ldots,t_{j}\right)  $ as an increasing
list of the set $S\setminus\left\{  m\right\}  $. In other words, $\left(
t_{1},t_{2},\ldots,t_{j}\right)  $ is a list of elements of $S\setminus
\left\{  m\right\}  $ such that $S\setminus\left\{  m\right\}  =\left\{
t_{1},t_{2},\ldots,t_{j}\right\}  $ and $t_{1}<t_{2}<\cdots<t_{j}$ (by the
definition of an \textquotedblleft increasing list\textquotedblright).

We claim that%
\begin{equation}
t_{1}<t_{2}<\cdots<t_{j+1}. \label{pf.thm.ind.inclist.unex.1}%
\end{equation}

[\textit{Proof of (\ref{pf.thm.ind.inclist.unex.1}):} If $j+1\leq1$, then the
chain of inequalities (\ref{pf.thm.ind.inclist.unex.1}) is vacuously true
(since it contains no inequality signs). Thus, for the rest of this proof of
(\ref{pf.thm.ind.inclist.unex.1}), we WLOG assume that we don't have
$j+1\leq1$. Hence, $j+1>1$, so that $j>0$ and thus $j\geq1$ (since $j$ is an
integer). Hence, $t_{j}$ is well-defined. We have $j\in\left\{  1,2,\ldots
,j\right\}  $ (since $j\geq1$) and thus $t_{j}\in\left\{  t_{1},t_{2}%
,\ldots,t_{j}\right\}  =S\setminus\left\{  m\right\}  \subseteq S$. Hence,
(\ref{eq.ind.max.def-max.2}) (applied to $t=t_{j}$) yields $\max S\geq t_{j}$.
Hence, $t_{j}\leq\max S=m$. Moreover, $t_{j}\notin\left\{  m\right\}  $ (since
$t_{j}\in S\setminus\left\{  m\right\}  $); in other words, $t_{j}\neq m$.
Combining this with $t_{j}\leq m$, we obtain $t_{j}<m=t_{j+1}$. Combining the
chain of inequalities $t_{1}<t_{2}<\cdots<t_{j}$ with the single inequality
$t_{j}<t_{j+1}$, we obtain the longer chain of inequalities $t_{1}%
<t_{2}<\cdots<t_{j}<t_{j+1}$. In other words, $t_{1}<t_{2}<\cdots<t_{j+1}$.
This proves (\ref{pf.thm.ind.inclist.unex.1}).]

Next, we shall prove that%
\begin{equation}
S=\left\{  t_{1},t_{2},\ldots,t_{j+1}\right\}  .
\label{pf.thm.ind.inclist.unex.2}%
\end{equation}

[\textit{Proof of (\ref{pf.thm.ind.inclist.unex.2}):} We have $\left(
S\setminus\left\{  m\right\}  \right)  \cup\left\{  m\right\}  =S\cup\left\{
m\right\}  =S$ (since $m\in S$). Thus,%
\begin{align*}
S  &  =\left(  \underbrace{S\setminus\left\{  m\right\}  }_{=\left\{
t_{1},t_{2},\ldots,t_{j}\right\}  }\right)  \cup\left\{  \underbrace{m}%
_{=t_{j+1}}\right\}  =\left\{  t_{1},t_{2},\ldots,t_{j}\right\}  \cup\left\{
t_{j+1}\right\}  =\left\{  t_{1},t_{2},\ldots,t_{j},t_{j+1}\right\} \\
&  =\left\{  t_{1},t_{2},\ldots,t_{j+1}\right\}  .
\end{align*}
This proves (\ref{pf.thm.ind.inclist.unex.2}).]

Clearly, $t_{1},t_{2},\ldots,t_{j+1}$ are elements of the set $\left\{
t_{1},t_{2},\ldots,t_{j+1}\right\}  $. In other words, $t_{1},t_{2}%
,\ldots,t_{j+1}$ are elements of the set $S$ (since $S=\left\{  t_{1}%
,t_{2},\ldots,t_{j+1}\right\}  $).

Hence, $\left(  t_{1},t_{2},\ldots,t_{j+1}\right)  $ is a list of elements of
$S$. Thus, $\left(  t_{1},t_{2},\ldots,t_{j+1}\right)  $ is a list of elements
of $S$ such that $S=\left\{  t_{1},t_{2},\ldots,t_{j+1}\right\}  $ (by
(\ref{pf.thm.ind.inclist.unex.2})) and $t_{1}<t_{2}<\cdots<t_{j+1}$ (by
(\ref{pf.thm.ind.inclist.unex.1})). In other words, $\left(  t_{1}%
,t_{2},\ldots,t_{j+1}\right)  $ is an increasing list of $S$ (by the
definition of an \textquotedblleft increasing list\textquotedblright). Hence,
the set $S$ has \textbf{at least} one increasing list (namely, $\left(
t_{1},t_{2},\ldots,t_{j+1}\right)  $).

We shall next show that $\left(  t_{1},t_{2},\ldots,t_{j+1}\right)  $ is the
only increasing list of $S$. Indeed, let $\left(  s_{1},s_{2},\ldots
,s_{k}\right)  $ be any increasing list of $S$. Then, Proposition
\ref{prop.ind.inclist.nonempty1} \textbf{(a)} shows that $k\geq1$ and
$s_{k}=m$. Also, Proposition \ref{prop.ind.inclist.nonempty1} \textbf{(b)}
shows that the list $\left(  s_{1},s_{2},\ldots,s_{k-1}\right)  $ is an
increasing list of $S\setminus\left\{  m\right\}  $.

But recall that $S\setminus\left\{  m\right\}  $ has exactly one increasing
list. Thus, in particular, $S\setminus\left\{  m\right\}  $ has \textbf{at
most} one increasing list. In other words, any two increasing lists of
$S\setminus\left\{  m\right\}  $ are equal. Hence, the lists $\left(
s_{1},s_{2},\ldots,s_{k-1}\right)  $ and $\left(  t_{1},t_{2},\ldots
,t_{j}\right)  $ must be equal (since both of these lists are increasing lists
of $S\setminus\left\{  m\right\}  $). In other words, $\left(  s_{1}%
,s_{2},\ldots,s_{k-1}\right)  =\left(  t_{1},t_{2},\ldots,t_{j}\right)  $. In
other words, $k-1=j$ and%
\begin{equation}
\left(  s_{i}=t_{i}\text{ for each }i\in\left\{  1,2,\ldots,k-1\right\}
\right)  . \label{pf.thm.ind.inclist.unex.3}%
\end{equation}

From $k-1=j$, we obtain $k=j+1$. Hence, $t_{k}=t_{j+1}=m$. Next, we claim that%
\begin{equation}
s_{i}=t_{i}\text{ for each }i\in\left\{  1,2,\ldots,k\right\}  .
\label{pf.thm.ind.inclist.unex.4}%
\end{equation}

[\textit{Proof of (\ref{pf.thm.ind.inclist.unex.4}):} Let $i\in\left\{
1,2,\ldots,k\right\}  $. We must prove that $s_{i}=t_{i}$. If $i\in\left\{
1,2,\ldots,k-1\right\}  $, then this follows from
(\ref{pf.thm.ind.inclist.unex.3}). Hence, for the rest of this proof, we WLOG
assume that we don't have $i\in\left\{  1,2,\ldots,k-1\right\}  $. Hence,
$i\notin\left\{  1,2,\ldots,k-1\right\}  $. Combining $i\in\left\{
1,2,\ldots,k\right\}  $ with $i\notin\left\{  1,2,\ldots,k-1\right\}  $, we
obtain%
\[
i\in\left\{  1,2,\ldots,k\right\}  \setminus\left\{  1,2,\ldots,k-1\right\}
=\left\{  k\right\}  .
\]
In other words, $i=k$. Hence, $s_{i}=s_{k}=m=t_{k}$ (since $t_{k}=m$). In view
of $k=i$, this rewrites as $s_{i}=t_{i}$. This proves
(\ref{pf.thm.ind.inclist.unex.4}).]

From (\ref{pf.thm.ind.inclist.unex.4}), we obtain $\left(  s_{1},s_{2}%
,\ldots,s_{k}\right)  =\left(  t_{1},t_{2},\ldots,t_{k}\right)  =\left(
t_{1},t_{2},\ldots,t_{j+1}\right)  $ (since $k=j+1$).

Now, forget that we fixed $\left(  s_{1},s_{2},\ldots,s_{k}\right)  $. We thus
have proven that if $\left(  s_{1},s_{2},\ldots,s_{k}\right)  $ is any
increasing list of $S$, then $\left(  s_{1},s_{2},\ldots,s_{k}\right)
=\left(  t_{1},t_{2},\ldots,t_{j+1}\right)  $. In other words, any increasing
list of $S$ equals $\left(  t_{1},t_{2},\ldots,t_{j+1}\right)  $. Thus, the
set $S$ has \textbf{at most} one increasing list. Since we also know that the
set $S$ has \textbf{at least} one increasing list, we thus conclude that $S$
has exactly one increasing list.

Now, forget that we fixed $S$. We thus have shown that%
\[
\left(
\begin{array}
[c]{c}%
\text{if }S\text{ is a finite set of integers satisfying }\left\vert
S\right\vert =g+1\text{,}\\
\text{then }S\text{ has exactly one increasing list}%
\end{array}
\right)  .
\]
In other words, Theorem \ref{thm.ind.inclist.unex} holds under the condition
that $\left\vert S\right\vert =g+1$. This completes the induction step. Hence,
Theorem \ref{thm.ind.inclist.unex} is proven by induction.
\end{proof}

\begin{definition}
\label{def.ind.inclist}Let $S$ be a finite set of integers. Theorem
\ref{thm.ind.inclist.unex} shows that $S$ has exactly one increasing list.
This increasing list is called \textit{the increasing list} of $S$. It is also
called \textit{the list of all elements of }$S$\textit{ in increasing order
(with no repetitions)}. (The latter name, of course, is descriptive.)

The increasing list of $S$ has length $\left\vert S\right\vert $. (Indeed, if
we denote this increasing list by $\left(  s_{1},s_{2},\ldots,s_{k}\right)  $,
then its length is $k=\left\vert S\right\vert $, because Proposition
\ref{prop.ind.inclist.size} \textbf{(b)} shows that $\left\vert S\right\vert
=k$.)

For each $j\in\left\{  1,2,\ldots,\left\vert S\right\vert \right\}  $, we
define the $j$\textit{-th smallest element of }$S$ to be the $j$-th entry of
the increasing list of $S$. In other words, if $\left(  s_{1},s_{2}%
,\ldots,s_{k}\right)  $ is the increasing list of $S$, then the $j$-th
smallest element of $S$ is $s_{j}$. Some say \textquotedblleft$j$-th lowest
element of $S$\textquotedblright\ instead of \textquotedblleft$j$-th smallest
element of $S$\textquotedblright.
\end{definition}

\begin{remark}
\textbf{(a)} Clearly, we can replace the word \textquotedblleft
integer\textquotedblright\ by \textquotedblleft rational
number\textquotedblright\ or by \textquotedblleft real
number\textquotedblright\ in Proposition \ref{prop.mod.chain-ineq}, Corollary
\ref{cor.mod.chain-ineq2}, Corollary \ref{cor.mod.chain-ineq3}, Definition
\ref{def.ind.inclist0}, Theorem \ref{thm.ind.inclist.unex}, Proposition
\ref{prop.ind.inclist.size}, Proposition \ref{prop.ind.inclist.empty},
Proposition \ref{prop.ind.inclist.nonempty1} and Definition
\ref{def.ind.inclist}, because we have not used any properties specific to integers.

\textbf{(b)} If we replace all the \textquotedblleft$<$\textquotedblright%
\ signs in Definition \ref{def.ind.inclist0} by \textquotedblleft%
$>$\textquotedblright\ signs, then we obtain the notion of a
\textit{decreasing list} of $S$. There are straightforward analogues of
Theorem \ref{thm.ind.inclist.unex}, Proposition \ref{prop.ind.inclist.size},
Proposition \ref{prop.ind.inclist.empty} and Proposition
\ref{prop.ind.inclist.nonempty1} for decreasing lists (where, of course, the
analogue of Proposition \ref{prop.ind.inclist.nonempty1} uses $\min S$ instead
of $\max S$). Thus, we can state an analogue of Definition
\ref{def.ind.inclist} as well. In this analogue, the word \textquotedblleft
increasing\textquotedblright\ is replaced by \textquotedblleft
decreasing\textquotedblright\ everywhere, the word \textquotedblleft
smallest\textquotedblright\ is replaced by \textquotedblleft
largest\textquotedblright, and the word \textquotedblleft
lowest\textquotedblright\ is replaced by \textquotedblleft
highest\textquotedblright.

\textbf{(c)} That said, the decreasing list and the increasing list are
closely related: If $S$ is a finite set of integers (or rational numbers, or
real numbers), and if $\left(  s_{1},s_{2},\ldots,s_{k}\right)  $ is the
increasing list of $S$, then $\left(  s_{k},s_{k-1},\ldots,s_{1}\right)  $ is
the decreasing list of $S$. (The proof is very simple.)

\textbf{(d)} Let $S$ be a nonempty finite set of integers (or rational
numbers, or real numbers), and let $\left(  s_{1},s_{2},\ldots,s_{k}\right)  $
be the increasing list of $S$. Proposition \ref{prop.ind.inclist.nonempty1}
\textbf{(a)} (applied to $m=\max S$) shows that $k\geq1$ and $s_{k}=\max S$. A
similar argument can be used to show that $s_{1}=\min S$. Thus, the increasing
list of $S$ begins with the smallest element of $S$ and ends with the largest
element of $S$ (as one would expect).
\end{remark}

\subsection{Induction with shifted base}

\subsubsection{Induction starting at $g$}

All the induction proofs we have done so far were applications of Theorem
\ref{thm.ind.IP0} (even though we have often written them up in ways that hide
the exact statements $\mathcal{A}\left(  n\right)  $ to which Theorem
\ref{thm.ind.IP0} is being applied). We are soon going to see several other
\textquotedblleft induction principles\textquotedblright\ which can also be
used to make proofs. Unlike Theorem \ref{thm.ind.IP0}, these other principles
need not be taken on trust; instead, they can themselves be proven using
Theorem \ref{thm.ind.IP0}. Thus, they merely offer convenience, not new
logical opportunities.

Our first such \textquotedblleft alternative induction
principle\textquotedblright\ is Theorem \ref{thm.ind.IPg} below. First, we
introduce a simple notation:

\begin{definition}
Let $g\in\mathbb{Z}$. Then, $\mathbb{Z}_{\geq g}$ denotes the set $\left\{
g,g+1,g+2,\ldots\right\}  $; this is the set of all integers that are $\geq g$.
\end{definition}

For example, $\mathbb{Z}_{\geq0}=\left\{  0,1,2,\ldots\right\}  =\mathbb{N}$
is the set of all nonnegative integers, whereas $\mathbb{Z}_{\geq1}=\left\{
1,2,3,\ldots\right\}  $ is the set of all positive integers.

Now, we state our first \textquotedblleft alternative induction
principle\textquotedblright:

\begin{theorem}
\label{thm.ind.IPg}Let $g\in\mathbb{Z}$. For each $n\in\mathbb{Z}_{\geq g}$,
let $\mathcal{A}\left(  n\right)  $ be a logical statement.

Assume the following:

\begin{statement}
\textit{Assumption 1:} The statement $\mathcal{A}\left(  g\right)  $ holds.
\end{statement}

\begin{statement}
\textit{Assumption 2:} If $m\in\mathbb{Z}_{\geq g}$ is such that
$\mathcal{A}\left(  m\right)  $ holds, then $\mathcal{A}\left(  m+1\right)  $
also holds.
\end{statement}

Then, $\mathcal{A}\left(  n\right)  $ holds for each $n\in\mathbb{Z}_{\geq g}$.
\end{theorem}

Again, Theorem \ref{thm.ind.IPg} is intuitively clear: For example, if you
have $g=4$, and you want to prove (under the assumptions of Theorem
\ref{thm.ind.IPg}) that $\mathcal{A}\left(  8\right)  $ holds, you can argue
as follows:

\begin{itemize}
\item By Assumption 1, the statement $\mathcal{A}\left(  4\right)  $ holds.

\item Thus, by Assumption 2 (applied to $m=4$), the statement $\mathcal{A}%
\left(  5\right)  $ holds.

\item Thus, by Assumption 2 (applied to $m=5$), the statement $\mathcal{A}%
\left(  6\right)  $ holds.

\item Thus, by Assumption 2 (applied to $m=6$), the statement $\mathcal{A}%
\left(  7\right)  $ holds.

\item Thus, by Assumption 2 (applied to $m=7$), the statement $\mathcal{A}%
\left(  8\right)  $ holds.
\end{itemize}

A similar (but longer) argument shows that the statement $\mathcal{A}\left(
9\right)  $ holds; likewise, $\mathcal{A}\left(  n\right)  $ can be shown to
hold for each $n\in\mathbb{Z}_{\geq g}$ by means of an argument that takes
$n-g+1$ steps.

Theorem \ref{thm.ind.IPg} generalizes Theorem \ref{thm.ind.IP0}. Indeed,
Theorem \ref{thm.ind.IP0} is the particular case of Theorem \ref{thm.ind.IPg}
for $g=0$ (since $\mathbb{Z}_{\geq0}=\mathbb{N}$). However, Theorem
\ref{thm.ind.IPg} can also be derived from Theorem \ref{thm.ind.IP0}. In order
to do this, we essentially need to \textquotedblleft shift\textquotedblright%
\ the index $n$ in Theorem \ref{thm.ind.IPg} down by $g$ -- that is, we need
to rename our sequence $\left(  \mathcal{A}\left(  g\right)  ,\mathcal{A}%
\left(  g+1\right)  ,\mathcal{A}\left(  g+2\right)  ,\ldots\right)  $ of
statements as $\left(  \mathcal{B}\left(  0\right)  ,\mathcal{B}\left(
1\right)  ,\mathcal{B}\left(  2\right)  ,\ldots\right)  $, and apply Theorem
\ref{thm.ind.IP0} to $\mathcal{B}\left(  n\right)  $ instead of $\mathcal{A}%
\left(  n\right)  $. In order to make this renaming procedure rigorous, let us
first restate Theorem \ref{thm.ind.IP0} as follows:

\begin{corollary}
\label{cor.ind.IP0.renamed}For each $n\in\mathbb{N}$, let $\mathcal{B}\left(
n\right)  $ be a logical statement.

Assume the following:

\begin{statement}
\textit{Assumption A:} The statement $\mathcal{B}\left(  0\right)  $ holds.
\end{statement}

\begin{statement}
\textit{Assumption B:} If $p\in\mathbb{N}$ is such that $\mathcal{B}\left(
p\right)  $ holds, then $\mathcal{B}\left(  p+1\right)  $ also holds.
\end{statement}

Then, $\mathcal{B}\left(  n\right)  $ holds for each $n\in\mathbb{N}$.
\end{corollary}

\begin{proof}
[Proof of Corollary \ref{cor.ind.IP0.renamed}.]Corollary
\ref{cor.ind.IP0.renamed} is exactly Theorem \ref{thm.ind.IP0}, except that
some names have been changed:

\begin{itemize}
\item The statements $\mathcal{A}\left(  n\right)  $ have been renamed as
$\mathcal{B}\left(  n\right)  $.

\item Assumption 1 and Assumption 2 have been renamed as Assumption A and
Assumption B.

\item The variable $m$ in Assumption B has been renamed as $p$.
\end{itemize}

Thus, Corollary \ref{cor.ind.IP0.renamed} holds (since Theorem
\ref{thm.ind.IP0} holds).
\end{proof}

Let us now derive Theorem \ref{thm.ind.IPg} from Theorem \ref{thm.ind.IP0}:

\begin{proof}
[Proof of Theorem \ref{thm.ind.IPg}.]For any $n\in\mathbb{N}$, we have
$n+g\in\mathbb{Z}_{\geq g}$\ \ \ \ \footnote{\textit{Proof.} Let
$n\in\mathbb{N}$. Thus, $n\geq0$, so that $\underbrace{n}_{\geq0}+g\geq0+g=g$.
Hence, $n+g$ is an integer $\geq g$. In other words, $n+g\in\mathbb{Z}_{\geq
g}$ (since $\mathbb{Z}_{\geq g}$ is the set of all integers that are $\geq
g$). Qed.}. Hence, for each $n\in\mathbb{N}$, we can define a logical
statement $\mathcal{B}\left(  n\right)  $ by%
\[
\mathcal{B}\left(  n\right)  =\mathcal{A}\left(  n+g\right)  .
\]
Consider this $\mathcal{B}\left(  n\right)  $.

Now, let us consider the Assumptions A and B from Corollary
\ref{cor.ind.IP0.renamed}. We claim that both of these assumptions are satisfied.

Indeed, the statement $\mathcal{A}\left(  g\right)  $ holds (by Assumption 1).
But the definition of the statement $\mathcal{B}\left(  0\right)  $ shows that
$\mathcal{B}\left(  0\right)  =\mathcal{A}\left(  0+g\right)  =\mathcal{A}%
\left(  g\right)  $. Hence, the statement $\mathcal{B}\left(  0\right)  $
holds (since the statement $\mathcal{A}\left(  g\right)  $ holds). In other
words, Assumption A is satisfied.

Now, we shall show that Assumption B is satisfied. Indeed, let $p\in
\mathbb{N}$ be such that $\mathcal{B}\left(  p\right)  $ holds. The definition
of the statement $\mathcal{B}\left(  p\right)  $ shows that $\mathcal{B}%
\left(  p\right)  =\mathcal{A}\left(  p+g\right)  $. Hence, the statement
$\mathcal{A}\left(  p+g\right)  $ holds (since $\mathcal{B}\left(  p\right)  $ holds).

Also, $p\in\mathbb{N}$, so that $p\geq0$ and thus $p+g\geq g$. In other words,
$p+g\in\mathbb{Z}_{\geq g}$ (since $\mathbb{Z}_{\geq g}$ is the set of all
integers that are $\geq g$).

Recall that Assumption 2 holds. In other words, if $m\in\mathbb{Z}_{\geq g}$
is such that $\mathcal{A}\left(  m\right)  $ holds, then $\mathcal{A}\left(
m+1\right)  $ also holds. Applying this to $m=p+g$, we conclude that
$\mathcal{A}\left(  \left(  p+g\right)  +1\right)  $ holds (since
$\mathcal{A}\left(  p+g\right)  $ holds).

But the definition of $\mathcal{B}\left(  p+1\right)  $ yields $\mathcal{B}%
\left(  p+1\right)  =\mathcal{A}\left(  \underbrace{p+1+g}_{=\left(
p+g\right)  +1}\right)  =\mathcal{A}\left(  \left(  p+g\right)  +1\right)  $.
Hence, the statement $\mathcal{B}\left(  p+1\right)  $ holds (since the
statement $\mathcal{A}\left(  \left(  p+g\right)  +1\right)  $ holds).

Now, forget that we fixed $p$. We thus have shown that if $p\in\mathbb{N}$ is
such that $\mathcal{B}\left(  p\right)  $ holds, then $\mathcal{B}\left(
p+1\right)  $ also holds. In other words, Assumption B is satisfied.

We now know that both Assumption A and Assumption B are satisfied. Hence,
Corollary \ref{cor.ind.IP0.renamed} shows that
\begin{equation}
\mathcal{B}\left(  n\right)  \text{ holds for each }n\in\mathbb{N}.
\label{pf.thm.ind.IPg.at}%
\end{equation}

Now, let $n\in\mathbb{Z}_{\geq g}$. Thus, $n$ is an integer such that $n\geq
g$ (by the definition of $\mathbb{Z}_{\geq g}$). Hence, $n-g\geq0$, so that
$n-g\in\mathbb{N}$. Thus, (\ref{pf.thm.ind.IPg.at}) (applied to $n-g$ instead
of $n$) yields that $\mathcal{B}\left(  n-g\right)  $ holds. But the
definition of $\mathcal{B}\left(  n-g\right)  $ yields $\mathcal{B}\left(
n-g\right)  =\mathcal{A}\left(  \underbrace{\left(  n-g\right)  +g}%
_{=n}\right)  =\mathcal{A}\left(  n\right)  $. Hence, the statement
$\mathcal{A}\left(  n\right)  $ holds (since $\mathcal{B}\left(  n-g\right)  $ holds).

Now, forget that we fixed $n$. We thus have shown that $\mathcal{A}\left(
n\right)  $ holds for each $n\in\mathbb{Z}_{\geq g}$. This proves Theorem
\ref{thm.ind.IPg}.
\end{proof}

Theorem \ref{thm.ind.IPg} is called the \textit{principle of induction
starting at }$g$, and proofs that use it are usually called \textit{proofs by
induction} or \textit{induction proofs}. As with the standard induction
principle (Theorem \ref{thm.ind.IP0}), we don't usually explicitly cite
Theorem \ref{thm.ind.IPg}, but instead say certain words that signal that it
is being applied and that (ideally) also indicate what integer $g$ and what
statements $\mathcal{A}\left(  n\right)  $ it is being applied to\footnote{We
will explain this in Convention \ref{conv.ind.IPglang} below.}. However, for
our very first example of the use of Theorem \ref{thm.ind.IPg}, we are going
to reference it explicitly:

\begin{proposition}
\label{prop.mod.binom01}Let $a$ and $b$ be integers. Then, every positive
integer $n$ satisfies%
\begin{equation}
\left(  a+b\right)  ^{n}\equiv a^{n}+na^{n-1}b\operatorname{mod}b^{2}.
\label{eq.prop.mod.binom01.claim}%
\end{equation}

\end{proposition}

Note that we have chosen not to allow $n=0$ in Proposition
\ref{prop.mod.binom01}, because it is not clear what \textquotedblleft%
$a^{n-1}$\textquotedblright\ would mean when $n=0$ and $a=0$. (Recall that
$0^{0-1}=0^{-1}$ is not defined!) In truth, it is easy to convince oneself
that this is not a serious hindrance, since the expression \textquotedblleft%
$na^{n-1}$\textquotedblright\ has a meaningful interpretation even when its
sub-expression \textquotedblleft$a^{n-1}$\textquotedblright\ does not (one
just has to interpret it as $0$ when $n=0$, without regard to whether
\textquotedblleft$a^{n-1}$\textquotedblright\ is well-defined). Nevertheless,
we prefer to rule out the case of $n=0$ by requiring $n$ to be positive, in
order to avoid having to discuss such questions of interpretation. (Of course,
this also gives us an excuse to apply Theorem \ref{thm.ind.IPg} instead of the
old Theorem \ref{thm.ind.IP0}.)

\begin{proof}
[Proof of Proposition \ref{prop.mod.binom01}.]For each $n\in\mathbb{Z}_{\geq
1}$, we let $\mathcal{A}\left(  n\right)  $ be the statement%
\[
\left(  \left(  a+b\right)  ^{n}\equiv a^{n}+na^{n-1}b\operatorname{mod}%
b^{2}\right)  .
\]
Our next goal is to prove the statement $\mathcal{A}\left(  n\right)  $ for
each $n\in\mathbb{Z}_{\geq1}$.

We first notice that the statement $\mathcal{A}\left(  1\right)  $
holds\footnote{\textit{Proof.} We have $\left(  a+b\right)  ^{1}=a+b$.
Comparing this with $\underbrace{a^{1}}_{=a}+1\underbrace{a^{1-1}}_{=a^{0}%
=1}b=a+b$, we obtain $\left(  a+b\right)  ^{1}=a^{1}+1a^{1-1}b$. Hence,
$\left(  a+b\right)  ^{1}\equiv a^{1}+1a^{1-1}b\operatorname{mod}b^{2}$. But
this is precisely the statement $\mathcal{A}\left(  1\right)  $ (since
$\mathcal{A}\left(  1\right)  $ is defined to be the statement $\left(
\left(  a+b\right)  ^{1}\equiv a^{1}+1a^{1-1}b\operatorname{mod}b^{2}\right)
$). Hence, the statement $\mathcal{A}\left(  1\right)  $ holds.}.

Now, we claim that
\begin{equation}
\text{if }m\in\mathbb{Z}_{\geq1}\text{ is such that }\mathcal{A}\left(
m\right)  \text{ holds, then }\mathcal{A}\left(  m+1\right)  \text{ also
holds.} \label{pf.prop.mod.binom01.step}%
\end{equation}

[\textit{Proof of (\ref{pf.prop.mod.binom01.step}):} Let $m\in\mathbb{Z}%
_{\geq1}$ be such that $\mathcal{A}\left(  m\right)  $ holds. We must show
that $\mathcal{A}\left(  m+1\right)  $ also holds.

We have assumed that $\mathcal{A}\left(  m\right)  $ holds. In other words,
\[
\left(  a+b\right)  ^{m}\equiv a^{m}+ma^{m-1}b\operatorname{mod}b^{2}%
\]
holds\footnote{because $\mathcal{A}\left(  m\right)  $ is defined to be the
statement $\left(  \left(  a+b\right)  ^{m}\equiv a^{m}+ma^{m-1}%
b\operatorname{mod}b^{2}\right)  $}. Now,%
\begin{align*}
\left(  a+b\right)  ^{m+1}  &  =\underbrace{\left(  a+b\right)  ^{m}}_{\equiv
a^{m}+ma^{m-1}b\operatorname{mod}b^{2}}\left(  a+b\right) \\
&  \equiv\left(  a^{m}+ma^{m-1}b\right)  \left(  a+b\right) \\
&  =\underbrace{a^{m}a}_{=a^{m+1}}+a^{m}b+m\underbrace{a^{m-1}ba}%
_{\substack{=a^{m-1}ab=a^{m}b\\\text{(since }a^{m-1}a=a^{m}\text{)}%
}}+\underbrace{ma^{m-1}bb}_{\substack{=ma^{m-1}b^{2}\equiv0\operatorname{mod}%
b^{2}\\\text{(since }b^{2}\mid ma^{m-1}b^{2}\text{)}}}\\
&  \equiv a^{m+1}+\underbrace{a^{m}b+ma^{m}b}_{=\left(  m+1\right)  a^{m}%
b}+0\\
&  =a^{m+1}+\left(  m+1\right)  \underbrace{a^{m}}_{\substack{=a^{\left(
m+1\right)  -1}\\\text{(since }m=\left(  m+1\right)  -1\text{)}}%
}b=a^{m+1}+\left(  m+1\right)  a^{\left(  m+1\right)  -1}b\operatorname{mod}%
b^{2}.
\end{align*}

So we have shown that $\left(  a+b\right)  ^{m+1}\equiv a^{m+1}+\left(
m+1\right)  a^{\left(  m+1\right)  -1}b\operatorname{mod}b^{2}$. But this is
precisely the statement $\mathcal{A}\left(  m+1\right)  $%
\ \ \ \ \footnote{because $\mathcal{A}\left(  m+1\right)  $ is defined to be
the statement $\left(  \left(  a+b\right)  ^{m+1}\equiv a^{m+1}+\left(
m+1\right)  a^{\left(  m+1\right)  -1}b\operatorname{mod}b^{2}\right)  $}.
Thus, the statement $\mathcal{A}\left(  m+1\right)  $ holds.

Now, forget that we fixed $m$. We thus have shown that if $m\in\mathbb{Z}%
_{\geq1}$ is such that $\mathcal{A}\left(  m\right)  $ holds, then
$\mathcal{A}\left(  m+1\right)  $ also holds. This proves
(\ref{pf.prop.mod.binom01.step}).]

Now, both assumptions of Theorem \ref{thm.ind.IPg} (applied to $g=1$) are
satisfied (indeed, Assumption 1 is satisfied because the statement
$\mathcal{A}\left(  1\right)  $ holds, whereas Assumption 2 is satisfied
because of (\ref{pf.prop.mod.binom01.step})). Thus, Theorem \ref{thm.ind.IPg}
(applied to $g=1$) shows that $\mathcal{A}\left(  n\right)  $ holds for each
$n\in\mathbb{Z}_{\geq1}$. In other words, $\left(  a+b\right)  ^{n}\equiv
a^{n}+na^{n-1}b\operatorname{mod}b^{2}$ holds for each $n\in\mathbb{Z}_{\geq
1}$ (since $\mathcal{A}\left(  n\right)  $ is the statement $\left(  \left(
a+b\right)  ^{n}\equiv a^{n}+na^{n-1}b\operatorname{mod}b^{2}\right)  $). In
other words, $\left(  a+b\right)  ^{n}\equiv a^{n}+na^{n-1}b\operatorname{mod}%
b^{2}$ holds for each positive integer $n$ (because the positive integers are
exactly the $n\in\mathbb{Z}_{\geq1}$). This proves Proposition
\ref{prop.mod.binom01}.
\end{proof}

\subsubsection{Conventions for writing proofs by induction starting at $g$}

Now, let us introduce some standard language that is commonly used in proofs
by induction starting at $g$:

\begin{convention}
\label{conv.ind.IPglang}Let $g\in\mathbb{Z}$. For each $n\in\mathbb{Z}_{\geq
g}$, let $\mathcal{A}\left(  n\right)  $ be a logical statement. Assume that
you want to prove that $\mathcal{A}\left(  n\right)  $ holds for each
$n\in\mathbb{Z}_{\geq g}$.

Theorem \ref{thm.ind.IPg} offers the following strategy for proving this:
First show that Assumption 1 of Theorem \ref{thm.ind.IPg} is satisfied; then,
show that Assumption 2 of Theorem \ref{thm.ind.IPg} is satisfied; then,
Theorem \ref{thm.ind.IPg} automatically completes your proof.

A proof that follows this strategy is called a \textit{proof by induction on
}$n$ (or \textit{proof by induction over }$n$) \textit{starting at }$g$ or
(less precisely) an \textit{inductive proof}. Most of the time, the words
\textquotedblleft starting at $g$\textquotedblright\ are omitted, since they
merely repeat what is clear from the context anyway: For example, if you make
a claim about all integers $n\geq3$, and you say that you are proving it by
induction on $n$, then it is clear that you are using induction on $n$
starting at $3$. (And if this isn't clear from the claim, then the induction
base will make it clear.)

The proof that Assumption 1 is satisfied is called the \textit{induction base}
(or \textit{base case}) of the proof. The proof that Assumption 2 is satisfied
is called the \textit{induction step} of the proof.

In order to prove that Assumption 2 is satisfied, you will usually want to fix
an $m\in\mathbb{Z}_{\geq g}$ such that $\mathcal{A}\left(  m\right)  $ holds,
and then prove that $\mathcal{A}\left(  m+1\right)  $ holds. In other words,
you will usually want to fix $m\in\mathbb{Z}_{\geq g}$, assume that
$\mathcal{A}\left(  m\right)  $ holds, and then prove that $\mathcal{A}\left(
m+1\right)  $ holds. When doing so, it is common to refer to the assumption
that $\mathcal{A}\left(  m\right)  $ holds as the \textit{induction
hypothesis} (or \textit{induction assumption}).
\end{convention}

Unsurprisingly, this language parallels the language introduced in Convention
\ref{conv.ind.IP0lang} for proofs by \textquotedblleft
standard\textquotedblright\ induction.

Again, we can shorten our inductive proofs by omitting some sentences that
convey no information. In particular, we can leave out the explicit definition
of the statement $\mathcal{A}\left(  n\right)  $ when this statement is
precisely the claim that we are proving (without the \textquotedblleft for
each $n\in\mathbb{Z}_{\geq g}$\textquotedblright\ part). Thus, we can rewrite
our above proof of Proposition \ref{prop.mod.binom01} as follows:

\begin{proof}
[Proof of Proposition \ref{prop.mod.binom01} (second version).]We must prove
(\ref{eq.prop.mod.binom01.claim}) for every positive integer $n$. In other
words, we must prove (\ref{eq.prop.mod.binom01.claim}) for every
$n\in\mathbb{Z}_{\geq1}$ (since the positive integers are precisely the
$n\in\mathbb{Z}_{\geq1}$). We shall prove this by induction on $n$ starting at
$1$:

\textit{Induction base:} We have $\left(  a+b\right)  ^{1}=a+b$. Comparing
this with $\underbrace{a^{1}}_{=a}+1\underbrace{a^{1-1}}_{=a^{0}=1}b=a+b$, we
obtain $\left(  a+b\right)  ^{1}=a^{1}+1a^{1-1}b$. Hence, $\left(  a+b\right)
^{1}\equiv a^{1}+1a^{1-1}b\operatorname{mod}b^{2}$. In other words,
(\ref{eq.prop.mod.binom01.claim}) holds for $n=1$. This completes the
induction base.

\textit{Induction step:} Let $m\in\mathbb{Z}_{\geq1}$. Assume that
(\ref{eq.prop.mod.binom01.claim}) holds for $n=m$. We must show that
(\ref{eq.prop.mod.binom01.claim}) also holds for $n=m+1$.

We have assumed that (\ref{eq.prop.mod.binom01.claim}) holds for $n=m$. In
other words,
\[
\left(  a+b\right)  ^{m}\equiv a^{m}+ma^{m-1}b\operatorname{mod}b^{2}%
\]
holds. Now,%
\begin{align*}
\left(  a+b\right)  ^{m+1}  &  =\underbrace{\left(  a+b\right)  ^{m}}_{\equiv
a^{m}+ma^{m-1}b\operatorname{mod}b^{2}}\left(  a+b\right) \\
&  \equiv\left(  a^{m}+ma^{m-1}b\right)  \left(  a+b\right) \\
&  =\underbrace{a^{m}a}_{=a^{m+1}}+a^{m}b+m\underbrace{a^{m-1}ba}%
_{\substack{=a^{m-1}ab=a^{m}b\\\text{(since }a^{m-1}a=a^{m}\text{)}%
}}+\underbrace{ma^{m-1}bb}_{\substack{=ma^{m-1}b^{2}\equiv0\operatorname{mod}%
b^{2}\\\text{(since }b^{2}\mid ma^{m-1}b^{2}\text{)}}}\\
&  \equiv a^{m+1}+\underbrace{a^{m}b+ma^{m}b}_{=\left(  m+1\right)  a^{m}%
b}+0\\
&  =a^{m+1}+\left(  m+1\right)  \underbrace{a^{m}}_{\substack{=a^{\left(
m+1\right)  -1}\\\text{(since }m=\left(  m+1\right)  -1\text{)}}%
}b=a^{m+1}+\left(  m+1\right)  a^{\left(  m+1\right)  -1}b\operatorname{mod}%
b^{2}.
\end{align*}

So we have shown that $\left(  a+b\right)  ^{m+1}\equiv a^{m+1}+\left(
m+1\right)  a^{\left(  m+1\right)  -1}b\operatorname{mod}b^{2}$. In other
words, (\ref{eq.prop.mod.binom01.claim}) holds for $n=m+1$.

Now, forget that we fixed $m$. We thus have shown that if $m\in\mathbb{Z}%
_{\geq1}$ is such that (\ref{eq.prop.mod.binom01.claim}) holds for $n=m$, then
(\ref{eq.prop.mod.binom01.claim}) also holds for $n=m+1$. This completes the
induction step. Hence, (\ref{eq.prop.mod.binom01.claim}) is proven by
induction. This proves Proposition \ref{prop.mod.binom01}.
\end{proof}

Proposition \ref{prop.mod.binom01} can also be seen as a consequence of the
binomial formula (Proposition \ref{prop.binom.binomial} further below).

\subsubsection{More properties of congruences}

Let us use this occasion to show two corollaries of Proposition
\ref{prop.mod.binom01}:

\begin{corollary}
\label{cor.mod.lte1}Let $a$, $b$ and $n$ be three integers such that $a\equiv
b\operatorname{mod}n$. Let $d\in\mathbb{N}$ be such that $d\mid n$. Then,
$a^{d}\equiv b^{d}\operatorname{mod}nd$.
\end{corollary}

\begin{proof}
[Proof of Corollary \ref{cor.mod.lte1}.]We have $a\equiv b\operatorname{mod}%
n$. In other words, $a$ is congruent to $b$ modulo $n$. In other words, $n\mid
a-b$ (by the definition of \textquotedblleft congruent\textquotedblright). In
other words, there exists an integer $w$ such that $a-b=nw$. Consider this
$w$. From $a-b=nw$, we obtain $a=b+nw$. Also, $d\mid n$, thus $dn\mid nn$ (by
Proposition \ref{prop.div.acbc}, applied to $d$, $n$ and $n$ instead of $a$,
$b$ and $c$). On the other hand, $nn\mid\left(  nw\right)  ^{2}$ (since
$\left(  nw\right)  ^{2}=nwnw=nnww$). Hence, Proposition \ref{prop.div.trans}
(applied to $dn$, $nn$ and $\left(  nw\right)  ^{2}$ instead of $a$, $b$ and
$c$) yields $dn\mid\left(  nw\right)  ^{2}$ (since $dn\mid nn$ and
$nn\mid\left(  nw\right)  ^{2}$). In other words, $nd\mid\left(  nw\right)
^{2}$ (since $dn=nd$).

Next, we claim that%
\begin{equation}
nd\mid a^{d}-b^{d}. \label{pf.cor.mod.lte1.1}%
\end{equation}

[\textit{Proof of (\ref{pf.cor.mod.lte1.1}):} If $d=0$, then
(\ref{pf.cor.mod.lte1.1}) holds (because if $d=0$, then $a^{d}-b^{d}%
=\underbrace{a^{0}}_{=1}-\underbrace{b^{0}}_{=1}=1-1=0=0nd$, and thus $nd\mid
a^{d}-b^{d}$). Hence, for the rest of this proof of (\ref{pf.cor.mod.lte1.1}),
we WLOG assume that we don't have $d=0$. Thus, $d\neq0$. Hence, $d$ is a
positive integer (since $d\in\mathbb{N}$). Thus, Proposition
\ref{prop.mod.binom01} (applied to $d$, $b$ and $nw$ instead of $n$, $a$ and
$b$) yields
\[
\left(  b+nw\right)  ^{d}\equiv b^{d}+db^{d-1}nw\operatorname{mod}\left(
nw\right)  ^{2}.
\]
In view of $a=b+nw$, this rewrites as%
\[
a^{d}\equiv b^{d}+db^{d-1}nw\operatorname{mod}\left(  nw\right)  ^{2}.
\]
Hence, Proposition \ref{prop.mod.0} \textbf{(c)} (applied to $a^{d}$,
$b^{d}+db^{d-1}nw$, $\left(  nw\right)  ^{2}$ and $nd$ instead of $a$, $b$,
$n$ and $m$) yields
\[
a^{d}\equiv b^{d}+db^{d-1}nw\operatorname{mod}nd
\]
(since $nd\mid\left(  nw\right)  ^{2}$). Hence,%
\[
a^{d}\equiv b^{d}+\underbrace{db^{d-1}nw}_{\substack{=ndb^{d-1}w\equiv
0\operatorname{mod}nd\\\text{(since }nd\mid ndb^{d-1}w\text{)}}}\equiv
b^{d}+0=b^{d}\operatorname{mod}nd.
\]
In other words, $nd\mid a^{d}-b^{d}$. This proves (\ref{pf.cor.mod.lte1.1}).]

From (\ref{pf.cor.mod.lte1.1}), we immediately obtain $a^{d}\equiv
b^{d}\operatorname{mod}nd$ (by the definition of \textquotedblleft
congruent\textquotedblright). This proves Corollary \ref{cor.mod.lte1}.
\end{proof}

For the next corollary, we need a convention:

\begin{convention}
\label{conv.triple-power}Let $a$, $b$ and $c$ be three integers. Then, the
expression \textquotedblleft$a^{b^{c}}$\textquotedblright\ shall always be
interpreted as \textquotedblleft$a^{\left(  b^{c}\right)  }$\textquotedblright%
, never as \textquotedblleft$\left(  a^{b}\right)  ^{c}$\textquotedblright.
\end{convention}

Thus, for example, \textquotedblleft$3^{3^{3}}$\textquotedblright\ means
$3^{\left(  3^{3}\right)  }=3^{27}=\allowbreak7625\,\allowbreak597\,484\,987$,
not $\left(  3^{3}\right)  ^{3}=27^{3}=\allowbreak19\,683$. The reason for
this convention is that $\left(  a^{b}\right)  ^{c}$ can be simplified to
$a^{bc}$ and thus there is little use in having yet another notation for it.
Of course, this convention applies not only to integers, but to any other
numbers $a,b,c$.

We can now state the following fact, which is sometimes known as
\textquotedblleft lifting-the-exponent lemma\textquotedblright:

\begin{corollary}
\label{cor.mod.lte2}Let $n\in\mathbb{N}$. Let $a$ and $b$ be two integers such
that $a\equiv b\operatorname{mod}n$. Let $k\in\mathbb{N}$. Then,
\begin{equation}
a^{n^{k}}\equiv b^{n^{k}}\operatorname{mod}n^{k+1}.
\label{eq.cor.mod.lte2.claim}%
\end{equation}

\end{corollary}

We shall give two \textbf{different} proofs of Corollary \ref{cor.mod.lte2} by
induction on $k$, to illustrate once again the point (previously made in
Remark \ref{rmk.ind.abstract}) that we have a choice of what precise statement
we are proving by induction. In the first proof, the statement will be the
congruence (\ref{eq.cor.mod.lte2.claim}) for three \textbf{fixed} integers
$a$, $b$ and $n$, whereas in the second proof, it will be the statement%
\[
\left(  a^{n^{k}}\equiv b^{n^{k}}\operatorname{mod}n^{k+1}\text{ for
\textbf{all} integers }a\text{ and }b\text{ and \textbf{all} }n\in
\mathbb{N}\text{ satisfying }a\equiv b\operatorname{mod}n\right)  .
\]

\begin{proof}
[First proof of Corollary \ref{cor.mod.lte2}.]Forget that we fixed $k$. We
thus must prove (\ref{eq.cor.mod.lte2.claim}) for each $k\in\mathbb{N}$.

We shall prove this by induction on $k$:

\textit{Induction base:} We have $n^{0}=1$ and thus $a^{n^{0}}=a^{1}=a$.
Similarly, $b^{n^{0}}=b$. Thus, $a^{n^{0}}=a\equiv b=b^{n^{0}}%
\operatorname{mod}n$. In other words, $a^{n^{0}}\equiv b^{n^{0}}%
\operatorname{mod}n^{0+1}$ (since $n^{0+1}=n^{1}=n$). In other words,
(\ref{eq.cor.mod.lte2.claim}) holds for $k=0$. This completes the induction base.

\textit{Induction step:} Let $m\in\mathbb{N}$. Assume that
(\ref{eq.cor.mod.lte2.claim}) holds for $k=m$. We must prove that
(\ref{eq.cor.mod.lte2.claim}) holds for $k=m+1$.

We have $n^{m+1}=nn^{m}$. Hence, $n\mid n^{m+1}$.

We have assumed that (\ref{eq.cor.mod.lte2.claim}) holds for $k=m$. In other
words, we have%
\[
a^{n^{m}}\equiv b^{n^{m}}\operatorname{mod}n^{m+1}.
\]
Hence, Corollary \ref{cor.mod.lte1} (applied to $a^{n^{m}}$, $b^{n^{m}}$,
$n^{m+1}$ and $n$ instead of $a$, $b$, $n$ and $d$) yields%
\[
\left(  a^{n^{m}}\right)  ^{n}\equiv\left(  b^{n^{m}}\right)  ^{n}%
\operatorname{mod}n^{m+1}n.
\]
Now, $n^{m+1}=n^{m}n$, so that
\[
a^{n^{m+1}}=a^{n^{m}n}=\left(  a^{n^{m}}\right)  ^{n}\equiv\left(  b^{n^{m}%
}\right)  ^{n}=b^{n^{m}n}=b^{n^{m+1}}\operatorname{mod}n^{m+1}n
\]
(since $n^{m}n=n^{m+1}$). In view of $n^{m+1}n=n^{\left(  m+1\right)  +1}$,
this rewrites as%
\[
a^{n^{m+1}}\equiv b^{n^{m+1}}\operatorname{mod}n^{\left(  m+1\right)  +1}.
\]
In other words, (\ref{eq.cor.mod.lte2.claim}) holds for $k=m+1$. This
completes the induction step. Thus, (\ref{eq.cor.mod.lte2.claim}) is proven by
induction. Hence, Corollary \ref{cor.mod.lte2} holds.
\end{proof}

\begin{proof}
[Second proof of Corollary \ref{cor.mod.lte2}.]Forget that we fixed $a$, $b$,
$n$ and $k$. We thus must prove
\begin{equation}
\left(  a^{n^{k}}\equiv b^{n^{k}}\operatorname{mod}n^{k+1}\text{ for all
integers }a\text{ and }b\text{ and all }n\in\mathbb{N}\text{ satisfying
}a\equiv b\operatorname{mod}n\right)  \label{pf.cor.mod.lte2.pf2.goal}%
\end{equation}
for all $k\in\mathbb{N}$.

We shall prove this by induction on $k$:

\textit{Induction base:} Let $n\in\mathbb{N}$. Let $a$ and $b$ be two integers
such that $a\equiv b\operatorname{mod}n$. We have $n^{0}=1$ and thus
$a^{n^{0}}=a^{1}=a$. Similarly, $b^{n^{0}}=b$. Thus, $a^{n^{0}}=a\equiv
b=b^{n^{0}}\operatorname{mod}n$. In other words, $a^{n^{0}}\equiv b^{n^{0}%
}\operatorname{mod}n^{0+1}$ (since $n^{0+1}=n^{1}=n$).

Now, forget that we fixed $n$, $a$ and $b$. We thus have proven that
$a^{n^{0}}\equiv b^{n^{0}}\operatorname{mod}n^{0+1}$ for all integers $a$ and
$b$ and all $n\in\mathbb{N}$ satisfying $a\equiv b\operatorname{mod}n$. In
other words, (\ref{pf.cor.mod.lte2.pf2.goal}) holds for $k=0$. This completes
the induction base.

\textit{Induction step:} Let $m\in\mathbb{N}$. Assume that
(\ref{pf.cor.mod.lte2.pf2.goal}) holds for $k=m$. We must prove that
(\ref{pf.cor.mod.lte2.pf2.goal}) holds for $k=m+1$.

Let $n\in\mathbb{N}$. Let $a$ and $b$ be two integers such that $a\equiv
b\operatorname{mod}n$. Now,
\begin{align*}
\left(  n^{2}\right)  ^{m+1}  &  =n^{2\left(  m+1\right)  }=n^{\left(
m+2\right)  +m}\ \ \ \ \ \ \ \ \ \ \left(  \text{since }2\left(  m+1\right)
=\left(  m+2\right)  +m\right) \\
&  =n^{m+2}n^{m},
\end{align*}
so that $n^{m+2}\mid\left(  n^{2}\right)  ^{m+1}$.

We have $n\mid n$. Hence, Corollary \ref{cor.mod.lte1} (applied to $d=n$)
yields $a^{n}\equiv b^{n}\operatorname{mod}nn$. In other words, $a^{n}\equiv
b^{n}\operatorname{mod}n^{2}$ (since $nn=n^{2}$).

We have assumed that (\ref{pf.cor.mod.lte2.pf2.goal}) holds for $k=m$. Hence,
we can apply (\ref{pf.cor.mod.lte2.pf2.goal}) to $a^{n}$, $b^{n}$, $n^{2}$ and
$m$ instead of $a$, $b$, $n$ and $k$ (since $a^{n}\equiv b^{n}%
\operatorname{mod}n^{2}$). We thus conclude that%
\[
\left(  a^{n}\right)  ^{n^{m}}\equiv\left(  b^{n}\right)  ^{n^{m}%
}\operatorname{mod}\left(  n^{2}\right)  ^{m+1}.
\]
Now, $n^{m+1}=nn^{m}$, so that
\[
a^{n^{m+1}}=a^{nn^{m}}=\left(  a^{n}\right)  ^{n^{m}}\equiv\left(
b^{n}\right)  ^{n^{m}}=b^{nn^{m}}=b^{n^{m+1}}\operatorname{mod}\left(
n^{2}\right)  ^{m+1}%
\]
(since $nn^{m}=n^{m+1}$). Hence, Proposition \ref{prop.mod.0} \textbf{(c)}
(applied to $a^{n^{m+1}}$, $b^{n^{m+1}}$, $\left(  n^{2}\right)  ^{m+1}$ and
$n^{m+2}$ instead of $a$, $b$, $n$ and $m$) yields $a^{n^{m+1}}\equiv
b^{n^{m+1}}\operatorname{mod}n^{m+2}$ (since $n^{m+2}\mid\left(  n^{2}\right)
^{m+1}$). In view of $m+2=\left(  m+1\right)  +1$, this rewrites as%
\[
a^{n^{m+1}}\equiv b^{n^{m+1}}\operatorname{mod}n^{\left(  m+1\right)  +1}.
\]

Now, forget that we fixed $n$, $a$ and $b$. We thus have proven that
\newline$a^{n^{m+1}}\equiv b^{n^{m+1}}\operatorname{mod}n^{\left(  m+1\right)
+1}$ for all integers $a$ and $b$ and all $n\in\mathbb{N}$ satisfying $a\equiv
b\operatorname{mod}n$. In other words, (\ref{pf.cor.mod.lte2.pf2.goal}) holds
for $k=m+1$. This completes the induction step. Thus,
(\ref{pf.cor.mod.lte2.pf2.goal}) is proven by induction. Hence, Corollary
\ref{cor.mod.lte2} is proven again.
\end{proof}

\subsection{\label{sect.ind.SIP}Strong induction}

\subsubsection{The strong induction principle}

We shall now show another \textquotedblleft alternative induction
principle\textquotedblright, which is known as the \textit{strong induction
principle} because it feels stronger than Theorem \ref{thm.ind.IP0} (in the
sense that it appears to get the same conclusion from weaker assumptions).
Just as Theorem \ref{thm.ind.IPg}, this principle is not a new axiom, but
rather a consequence of the standard induction principle; we shall soon deduce
it from Theorem \ref{thm.ind.IPg}.

\begin{theorem}
\label{thm.ind.SIP}Let $g\in\mathbb{Z}$. For each $n\in\mathbb{Z}_{\geq g}$,
let $\mathcal{A}\left(  n\right)  $ be a logical statement.

Assume the following:

\begin{statement}
\textit{Assumption 1:} If $m\in\mathbb{Z}_{\geq g}$ is such that%
\[
\left(  \mathcal{A}\left(  n\right)  \text{ holds for every }n\in
\mathbb{Z}_{\geq g}\text{ satisfying }n<m\right)  ,
\]
then $\mathcal{A}\left(  m\right)  $ holds.
\end{statement}

Then, $\mathcal{A}\left(  n\right)  $ holds for each $n\in\mathbb{Z}_{\geq g}$.
\end{theorem}

Notice that Theorem \ref{thm.ind.SIP} has only one assumption (unlike Theorem
\ref{thm.ind.IP0} and Theorem \ref{thm.ind.IPg}). We shall soon see that this
one assumption \textquotedblleft incorporates\textquotedblright\ both an
induction base and an induction step.

Let us first explain why Theorem \ref{thm.ind.SIP} is intuitively clear. For
example, if you have $g=4$, and you want to prove (under the assumptions of
Theorem \ref{thm.ind.SIP}) that $\mathcal{A}\left(  7\right)  $ holds, you can
argue as follows:

\begin{itemize}
\item We know that $\mathcal{A}\left(  n\right)  $ holds for every
$n\in\mathbb{Z}_{\geq4}$ satisfying $n<4$. (Indeed, this is vacuously true,
since there is no $n\in\mathbb{Z}_{\geq4}$ satisfying $n<4$.)

Hence, Assumption 1 (applied to $m=4$) shows that the statement $\mathcal{A}%
\left(  4\right)  $ holds.

\item Thus, we know that $\mathcal{A}\left(  n\right)  $ holds for every
$n\in\mathbb{Z}_{\geq4}$ satisfying $n<5$ (because $\mathcal{A}\left(
4\right)  $ holds).

Hence, Assumption 1 (applied to $m=5$) shows that the statement $\mathcal{A}%
\left(  5\right)  $ holds.

\item Thus, we know that $\mathcal{A}\left(  n\right)  $ holds for every
$n\in\mathbb{Z}_{\geq4}$ satisfying $n<6$ (because $\mathcal{A}\left(
4\right)  $ and $\mathcal{A}\left(  5\right)  $ hold).

Hence, Assumption 1 (applied to $m=6$) shows that the statement $\mathcal{A}%
\left(  6\right)  $ holds.

\item Thus, we know that $\mathcal{A}\left(  n\right)  $ holds for every
$n\in\mathbb{Z}_{\geq4}$ satisfying $n<7$ (because $\mathcal{A}\left(
4\right)  $, $\mathcal{A}\left(  5\right)  $ and $\mathcal{A}\left(  6\right)
$ hold).

Hence, Assumption 1 (applied to $m=7$) shows that the statement $\mathcal{A}%
\left(  7\right)  $ holds.
\end{itemize}

A similar (but longer) argument shows that the statement $\mathcal{A}\left(
8\right)  $ holds; likewise, $\mathcal{A}\left(  n\right)  $ can be shown to
hold for each $n\in\mathbb{Z}_{\geq g}$ by means of an argument that takes
$n-g+1$ steps.

It is easy to see that Theorem \ref{thm.ind.SIP} generalizes Theorem
\ref{thm.ind.IPg} (because if the two Assumptions 1 and 2 of Theorem
\ref{thm.ind.IPg} hold, then so does Assumption 1 of Theorem \ref{thm.ind.SIP}%
). More interesting for us is the converse implication: We shall show that
Theorem \ref{thm.ind.SIP} can be derived from Theorem \ref{thm.ind.IPg}. This
will allow us to use Theorem \ref{thm.ind.SIP} without having to taking it on trust.

Before we derive Theorem \ref{thm.ind.SIP}, let us restate Theorem
\ref{thm.ind.IPg} as follows:

\begin{corollary}
\label{cor.ind.IPg.renamed}Let $g\in\mathbb{Z}$. For each $n\in\mathbb{Z}%
_{\geq g}$, let $\mathcal{B}\left(  n\right)  $ be a logical statement.

Assume the following:

\begin{statement}
\textit{Assumption A:} The statement $\mathcal{B}\left(  g\right)  $ holds.
\end{statement}

\begin{statement}
\textit{Assumption B:} If $p\in\mathbb{Z}_{\geq g}$ is such that
$\mathcal{B}\left(  p\right)  $ holds, then $\mathcal{B}\left(  p+1\right)  $
also holds.
\end{statement}

Then, $\mathcal{B}\left(  n\right)  $ holds for each $n\in\mathbb{Z}_{\geq g}$.
\end{corollary}

\begin{proof}
[Proof of Corollary \ref{cor.ind.IPg.renamed}.]Corollary
\ref{cor.ind.IPg.renamed} is exactly Theorem \ref{thm.ind.IPg}, except that
some names have been changed:

\begin{itemize}
\item The statements $\mathcal{A}\left(  n\right)  $ have been renamed as
$\mathcal{B}\left(  n\right)  $.

\item Assumption 1 and Assumption 2 have been renamed as Assumption A and
Assumption B.

\item The variable $m$ in Assumption B has been renamed as $p$.
\end{itemize}

Thus, Corollary \ref{cor.ind.IPg.renamed} holds (since Theorem
\ref{thm.ind.IPg} holds).
\end{proof}

Let us now derive Theorem \ref{thm.ind.SIP} from Theorem \ref{thm.ind.IPg}:

\begin{proof}
[Proof of Theorem \ref{thm.ind.SIP}.]For each $n\in\mathbb{Z}_{\geq g}$, we
let $\mathcal{B}\left(  n\right)  $ be the statement%
\[
\left(  \mathcal{A}\left(  q\right)  \text{ holds for every }q\in
\mathbb{Z}_{\geq g}\text{ satisfying }q<n\right)  .
\]

Now, let us consider the Assumptions A and B from Corollary
\ref{cor.ind.IPg.renamed}. We claim that both of these assumptions are satisfied.

The statement $\mathcal{B}\left(  g\right)  $ holds\footnote{\textit{Proof.}
Let $q\in\mathbb{Z}_{\geq g}$ be such that $q<g$. Then, $q\geq g$ (since
$q\in\mathbb{Z}_{\geq g}$); but this contradicts $q<g$.
\par
Now, forget that we fixed $q$. We thus have found a contradiction for each
$q\in\mathbb{Z}_{\geq g}$ satisfying $q<g$. Hence, there exists no
$q\in\mathbb{Z}_{\geq g}$ satisfying $q<g$. Thus, the statement
\[
\left(  \mathcal{A}\left(  q\right)  \text{ holds for every }q\in
\mathbb{Z}_{\geq g}\text{ satisfying }q<g\right)
\]
is vacuously true, and therefore true. In other words, the statement
$\mathcal{B}\left(  g\right)  $ is true (since $\mathcal{B}\left(  g\right)  $
is defined as the statement $\left(  \mathcal{A}\left(  q\right)  \text{ holds
for every }q\in\mathbb{Z}_{\geq g}\text{ satisfying }q<g\right)  $). Qed.}.
Thus, Assumption A is satisfied.

Next, let us prove that Assumption B is satisfied. Indeed, let $p\in
\mathbb{Z}_{\geq g}$ be such that $\mathcal{B}\left(  p\right)  $ holds. We
shall show that $\mathcal{B}\left(  p+1\right)  $ also holds.

Indeed, we have assumed that $\mathcal{B}\left(  p\right)  $ holds. In other
words,%
\begin{equation}
\mathcal{A}\left(  q\right)  \text{ holds for every }q\in\mathbb{Z}_{\geq
g}\text{ satisfying }q<p \label{pf.thm.ind.SIP.1}%
\end{equation}
(because the statement $\mathcal{B}\left(  p\right)  $ is defined as
\newline$\left(  \mathcal{A}\left(  q\right)  \text{ holds for every }%
q\in\mathbb{Z}_{\geq g}\text{ satisfying }q<p\right)  $). Renaming the
variable $q$ as $n$ in this statement, we conclude that%
\begin{equation}
\mathcal{A}\left(  n\right)  \text{ holds for every }n\in\mathbb{Z}_{\geq
g}\text{ satisfying }n<p. \label{pf.thm.ind.SIP.1n}%
\end{equation}
Hence, Assumption 1 (applied to $m=p$) yields that $\mathcal{A}\left(
p\right)  $ holds.

Now, we claim that%
\begin{equation}
\mathcal{A}\left(  q\right)  \text{ holds for every }q\in\mathbb{Z}_{\geq
g}\text{ satisfying }q<p+1. \label{pf.thm.ind.SIP.2}%
\end{equation}

[\textit{Proof of (\ref{pf.thm.ind.SIP.2}):} Let $q\in\mathbb{Z}_{\geq g}$ be
such that $q<p+1$. We must prove that $\mathcal{A}\left(  q\right)  $ holds.

If $q=p$, then this follows from the fact that $\mathcal{A}\left(  p\right)  $
holds. Hence, for the rest of this proof, we WLOG assume that we don't have
$q=p$. Thus, $q\neq p$. But $q<p+1$ and therefore $q\leq\left(  p+1\right)
-1$ (since $q$ and $p+1$ are integers). Hence, $q\leq\left(  p+1\right)
-1=p$. Combining this with $q\neq p$, we obtain $q<p$. Hence,
(\ref{pf.thm.ind.SIP.1}) shows that $\mathcal{A}\left(  q\right)  $ holds.
This completes the proof of (\ref{pf.thm.ind.SIP.2}).]

But the statement $\mathcal{B}\left(  p+1\right)  $ is defined as
\newline$\left(  \mathcal{A}\left(  q\right)  \text{ holds for every }%
q\in\mathbb{Z}_{\geq g}\text{ satisfying }q<p+1\right)  $. In other words, the
statement $\mathcal{B}\left(  p+1\right)  $ is precisely the statement
(\ref{pf.thm.ind.SIP.2}). Hence, the statement $\mathcal{B}\left(  p+1\right)
$ holds (since (\ref{pf.thm.ind.SIP.2}) holds).

Now, forget that we fixed $p$. We thus have shown that if $p\in\mathbb{Z}%
_{\geq g}$ is such that $\mathcal{B}\left(  p\right)  $ holds, then
$\mathcal{B}\left(  p+1\right)  $ also holds. In other words, Assumption B is satisfied.

We now know that both Assumption A and Assumption B are satisfied. Hence,
Corollary \ref{cor.ind.IPg.renamed} shows that
\begin{equation}
\mathcal{B}\left(  n\right)  \text{ holds for each }n\in\mathbb{Z}_{\geq g}.
\label{pf.thm.ind.SIP.at}%
\end{equation}

Now, let $n\in\mathbb{Z}_{\geq g}$. Thus, $n$ is an integer such that $n\geq
g$ (by the definition of $\mathbb{Z}_{\geq g}$). Hence, $n+1$ is also an
integer and satisfies $n+1\geq n\geq g$, so that $n+1\in\mathbb{Z}_{\geq g}$.
Hence, (\ref{pf.thm.ind.SIP.at}) (applied to $n+1$ instead of $n$) shows that
$\mathcal{B}\left(  n+1\right)  $ holds. In other words,%
\[
\mathcal{A}\left(  q\right)  \text{ holds for every }q\in\mathbb{Z}_{\geq
g}\text{ satisfying }q<n+1
\]
(because the statement $\mathcal{B}\left(  n+1\right)  $ is defined as
\newline$\left(  \mathcal{A}\left(  q\right)  \text{ holds for every }%
q\in\mathbb{Z}_{\geq g}\text{ satisfying }q<n+1\right)  $). We can apply this
to $q=n$ (because $n\in\mathbb{Z}_{\geq g}$ satisfies $n<n+1$), and conclude
that $\mathcal{A}\left(  n\right)  $ holds.

Now, forget that we fixed $n$. We thus have shown that $\mathcal{A}\left(
n\right)  $ holds for each $n\in\mathbb{Z}_{\geq g}$. This proves Theorem
\ref{thm.ind.SIP}.
\end{proof}

Thus, proving a sequence of statements $\mathcal{A}\left(  0\right)
,\mathcal{A}\left(  1\right)  ,\mathcal{A}\left(  2\right)  ,\ldots$ using
Theorem \ref{thm.ind.SIP} is tantamount to proving a slightly different
sequence of statements \newline$\mathcal{B}\left(  0\right)  ,\mathcal{B}%
\left(  1\right)  ,\mathcal{B}\left(  2\right)  ,\ldots$ using Corollary
\ref{cor.ind.IPg.renamed} and then deriving the former from the latter.

Theorem \ref{thm.ind.IPg} is called the \textit{principle of strong induction
starting at }$g$, and proofs that use it are usually called \textit{proofs by
strong induction}. We illustrate its use on the following easy property of the
Fibonacci sequence:

\begin{proposition}
\label{prop.ind.sip-fib}Let $\left(  f_{0},f_{1},f_{2},\ldots\right)  $ be the
Fibonacci sequence (defined as in Example \ref{exa.rec-seq.fib}). Then,%
\begin{equation}
f_{n}\leq2^{n-1} \label{eq.prop.ind.sip-fib.claim}%
\end{equation}
for each $n\in\mathbb{N}$.
\end{proposition}

\begin{proof}
[Proof of Proposition \ref{prop.ind.sip-fib}.]For each $n\in\mathbb{Z}_{\geq
0}$, we let $\mathcal{A}\left(  n\right)  $ be the statement $\left(
f_{n}\leq2^{n-1}\right)  $.

Thus, $\mathcal{A}\left(  0\right)  $ is the statement $\left(  f_{0}%
\leq2^{0-1}\right)  $; hence, this statement holds (since $f_{0}=0\leq2^{0-1}$).

Also, $\mathcal{A}\left(  1\right)  $ is the statement $\left(  f_{1}%
\leq2^{1-1}\right)  $ (by the definition of $\mathcal{A}\left(  1\right)  $);
hence, this statement also holds (since $f_{1}=1=2^{1-1}$).

Now, we claim the following:

\begin{statement}
\textit{Claim 1:} If $m\in\mathbb{Z}_{\geq0}$ is such that%
\[
\left(  \mathcal{A}\left(  n\right)  \text{ holds for every }n\in
\mathbb{Z}_{\geq0}\text{ satisfying }n<m\right)  ,
\]
then $\mathcal{A}\left(  m\right)  $ holds.
\end{statement}

[\textit{Proof of Claim 1:} Let $m\in\mathbb{Z}_{\geq0}$ be such that%
\begin{equation}
\left(  \mathcal{A}\left(  n\right)  \text{ holds for every }n\in
\mathbb{Z}_{\geq0}\text{ satisfying }n<m\right)  .
\label{pf.prop.ind.sip-fib.IH}%
\end{equation}
We must prove that $\mathcal{A}\left(  m\right)  $ holds.

\begin{vershort}
This is true if $m\in\left\{  0,1\right\}  $ (because we have shown that both
statements $\mathcal{A}\left(  0\right)  $ and $\mathcal{A}\left(  1\right)  $
hold). Thus, for the rest of the proof of Claim 1, we WLOG assume that we
don't have $m\in\left\{  0,1\right\}  $. Hence, $m\in\mathbb{N}\setminus
\left\{  0,1\right\}  =\left\{  2,3,4,\ldots\right\}  $, so that $m\geq2$.
\end{vershort}

\begin{verlong}
This is true if $m=0$ (because we have shown that the statement $\mathcal{A}%
\left(  0\right)  $ holds). Thus, for the rest of the proof of Claim 1, we
WLOG assume that we don't have $m=0$. Hence, $m\neq0$, so that $m\geq1$ (since
$m\in\mathbb{N}$).

We must prove that $\mathcal{A}\left(  m\right)  $ holds. This is true if
$m=1$ (because we have shown that the statement $\mathcal{A}\left(  1\right)
$ holds). Thus, for the rest of the proof of Claim 1, we WLOG assume that we
don't have $m=1$. Hence, $m\neq1$. Combining this with $m\geq1$, we obtain
$m>1$, so that $m\geq2$ (since $m$ is an integer).
\end{verlong}

From $m\geq2$, we conclude that $m-1\geq2-1=1\geq0$ and $m-2\geq2-2=0$. Thus,
both $m-1$ and $m-2$ belong to $\mathbb{N}$; therefore, $f_{m-1}$ and
$f_{m-2}$ are well-defined.

We have $m-1\in\mathbb{N}=\mathbb{Z}_{\geq0}$ and $m-1<m$. Hence,
(\ref{pf.prop.ind.sip-fib.IH}) (applied to $n=m-1$) yields that $\mathcal{A}%
\left(  m-1\right)  $ holds. In other words, $f_{m-1}\leq2^{\left(
m-1\right)  -1}$ (because this is what the statement $\mathcal{A}\left(
m-1\right)  $ says).

We have $m-2\in\mathbb{N}=\mathbb{Z}_{\geq0}$ and $m-2<m$. Hence,
(\ref{pf.prop.ind.sip-fib.IH}) (applied to $n=m-2$) yields that $\mathcal{A}%
\left(  m-2\right)  $ holds. In other words, $f_{m-2}\leq2^{\left(
m-2\right)  -1}$ (because this is what the statement $\mathcal{A}\left(
m-2\right)  $ says).

We have $\left(  m-1\right)  -1=m-2$ and thus $2^{\left(  m-1\right)
-1}=2^{m-2}=2\cdot2^{\left(  m-2\right)  -1}\geq2^{\left(  m-2\right)  -1}$
(since $2\cdot2^{\left(  m-2\right)  -1}-2^{\left(  m-2\right)  -1}=2^{\left(
m-2\right)  -1}\geq0$). Hence, $2^{\left(  m-2\right)  -1}\leq2^{\left(
m-1\right)  -1}$.

\begin{noncompile}
Old argument: Hence, $2^{\left(  m-2\right)  -1}\leq2^{\left(  m-1\right)
-1}$ (because if $a$, $b$ and $c$ are integers satisfying $a\leq b$ and
$c\geq1$, then $c^{a}\leq c^{b}$).
\end{noncompile}

But the recursive definition of the Fibonacci sequence yields $f_{m}%
=f_{m-1}+f_{m-2}$ (since $m\geq2$). Hence,%
\[
f_{m}=\underbrace{f_{m-1}}_{\leq2^{\left(  m-1\right)  -1}}%
+\underbrace{f_{m-2}}_{\leq2^{\left(  m-2\right)  -1}\leq2^{\left(
m-1\right)  -1}}\leq2^{\left(  m-1\right)  -1}+2^{\left(  m-1\right)
-1}=2\cdot2^{\left(  m-1\right)  -1}=2^{m-1}.
\]
In other words, the statement $\mathcal{A}\left(  m\right)  $ holds (since the
statement $\mathcal{A}\left(  m\right)  $ is defined to be $\left(  f_{m}%
\leq2^{m-1}\right)  $). This completes the proof of Claim 1.]

Claim 1 shows that Assumption 1 of Theorem \ref{thm.ind.SIP} (applied to
$g=0$) is satisfied. Hence, Theorem \ref{thm.ind.SIP} (applied to $g=0$) shows
that $\mathcal{A}\left(  n\right)  $ holds for each $n\in\mathbb{Z}_{\geq0}$.
In other words, $f_{n}\leq2^{n-1}$ holds for each $n\in\mathbb{Z}_{\geq0}$
(since the statement $\mathcal{A}\left(  n\right)  $ is defined to be $\left(
f_{n}\leq2^{n-1}\right)  $). In other words, $f_{n}\leq2^{n-1}$ holds for each
$n\in\mathbb{N}$ (since $\mathbb{Z}_{\geq0}=\mathbb{N}$). This proves
Proposition \ref{prop.ind.sip-fib}.
\end{proof}

\subsubsection{Conventions for writing strong induction proofs}

Again, when using the principle of strong induction, one commonly does not
directly cite Theorem \ref{thm.ind.SIP}; instead one uses the following language:

\begin{convention}
\label{conv.ind.SIPlang}Let $g\in\mathbb{Z}$. For each $n\in\mathbb{Z}_{\geq
g}$, let $\mathcal{A}\left(  n\right)  $ be a logical statement. Assume that
you want to prove that $\mathcal{A}\left(  n\right)  $ holds for each
$n\in\mathbb{Z}_{\geq g}$.

Theorem \ref{thm.ind.SIP} offers the following strategy for proving this: Show
that Assumption 1 of Theorem \ref{thm.ind.SIP} is satisfied; then, Theorem
\ref{thm.ind.SIP} automatically completes your proof.

A proof that follows this strategy is called a \textit{proof by strong
induction on }$n$ \textit{starting at }$g$. The proof that Assumption 1 is
satisfied is called the \textit{induction step} of the proof. This kind of
proof does not have an \textquotedblleft induction base\textquotedblright%
\ (unlike proofs that use Theorem \ref{thm.ind.IP0} or Theorem
\ref{thm.ind.IPg}).\footnotemark

In order to prove that Assumption 1 is satisfied, you will usually want to fix
an $m\in\mathbb{Z}_{\geq g}$ such that%
\begin{equation}
\left(  \mathcal{A}\left(  n\right)  \text{ holds for every }n\in
\mathbb{Z}_{\geq g}\text{ satisfying }n<m\right)  ,
\label{eq.conv.ind.SIPlang.IH}%
\end{equation}
and then prove that $\mathcal{A}\left(  m\right)  $ holds. In other words, you
will usually want to fix $m\in\mathbb{Z}_{\geq g}$, assume that
(\ref{eq.conv.ind.SIPlang.IH}) holds, and then prove that $\mathcal{A}\left(
m\right)  $ holds. When doing so, it is common to refer to the assumption that
(\ref{eq.conv.ind.SIPlang.IH}) holds as the \textit{induction hypothesis} (or
\textit{induction assumption}).
\end{convention}

\footnotetext{There is a version of strong induction which does include an
induction base (or even several). But the version we are using does not.}
Using this language, we can rewrite our above proof of Proposition
\ref{prop.ind.sip-fib} as follows:

\begin{proof}
[Proof of Proposition \ref{prop.ind.sip-fib} (second version).]For each
$n\in\mathbb{Z}_{\geq0}$, we let $\mathcal{A}\left(  n\right)  $ be the
statement $\left(  f_{n}\leq2^{n-1}\right)  $. Thus, our goal is to prove the
statement $\mathcal{A}\left(  n\right)  $ for each $n\in\mathbb{N}$. In other
words, our goal is to prove the statement $\mathcal{A}\left(  n\right)  $ for
each $n\in\mathbb{Z}_{\geq0}$ (since $\mathbb{N}=\mathbb{Z}_{\geq0}$).

We shall prove this by strong induction on $n$ starting at $0$:

\textit{Induction step:} Let $m\in\mathbb{Z}_{\geq0}$. Assume that%
\begin{equation}
\left(  \mathcal{A}\left(  n\right)  \text{ holds for every }n\in
\mathbb{Z}_{\geq0}\text{ satisfying }n<m\right)  .
\label{pf.prop.ind.sip-fib.2nd.IH}%
\end{equation}
We must then show that $\mathcal{A}\left(  m\right)  $ holds. In other words,
we must show that $f_{m}\leq2^{m-1}$ holds (since the statement $\mathcal{A}%
\left(  m\right)  $ is defined as $\left(  f_{m}\leq2^{m-1}\right)  $).

\begin{vershort}
This is true if $m=0$ (since $f_{0}=0\leq2^{0-1}$) and also true if $m=1$
(since $f_{1}=1=2^{1-1}$ and thus $f_{1}\leq2^{1-1}$). In other words, this is
true if $m\in\left\{  0,1\right\}  $. Thus, for the rest of the induction
step, we WLOG assume that we don't have $m\in\left\{  0,1\right\}  $. Hence,
$m\notin\left\{  0,1\right\}  $, so that $m\in\mathbb{N}\setminus\left\{
0,1\right\}  =\left\{  2,3,4,\ldots\right\}  $. Hence, $m\geq2$.
\end{vershort}

\begin{verlong}
This is true if $m=0$ (since $f_{0}=0\leq2^{0-1}$). Hence, for the rest of the
induction step, we WLOG assume that we don't have $m=0$. Hence, $m\neq0$, so
that $m\geq1$ (since $m\in\mathbb{Z}_{\geq0}$).

We must prove that $f_{m}\leq2^{m-1}$. This is true if $m=1$ (because
$f_{1}=1=2^{1-1}$ and thus $f_{1}\leq2^{1-1}$). Thus, for the rest of the
induction step, we WLOG assume that we don't have $m=1$. Hence, $m\neq1$.
Combining this with $m\geq1$, we obtain $m>1$, so that $m\geq2$ (since $m$ is
an integer).
\end{verlong}

From $m\geq2$, we conclude that $m-1\geq2-1=1\geq0$ and $m-2\geq2-2=0$. Thus,
both $m-1$ and $m-2$ belong to $\mathbb{N}$; therefore, $f_{m-1}$ and
$f_{m-2}$ are well-defined.

We have $m-1\in\mathbb{N}=\mathbb{Z}_{\geq0}$ and $m-1<m$. Hence,
(\ref{pf.prop.ind.sip-fib.2nd.IH}) (applied to $n=m-1$) yields that
$\mathcal{A}\left(  m-1\right)  $ holds. In other words, $f_{m-1}%
\leq2^{\left(  m-1\right)  -1}$ (because this is what the statement
$\mathcal{A}\left(  m-1\right)  $ says).

We have $m-2\in\mathbb{N}=\mathbb{Z}_{\geq0}$ and $m-2<m$. Hence,
(\ref{pf.prop.ind.sip-fib.2nd.IH}) (applied to $n=m-2$) yields that
$\mathcal{A}\left(  m-2\right)  $ holds. In other words, $f_{m-2}%
\leq2^{\left(  m-2\right)  -1}$ (because this is what the statement
$\mathcal{A}\left(  m-2\right)  $ says).

We have $\left(  m-1\right)  -1=m-2$ and thus $2^{\left(  m-1\right)
-1}=2^{m-2}=2\cdot2^{\left(  m-2\right)  -1}\geq2^{\left(  m-2\right)  -1}$
(since $2\cdot2^{\left(  m-2\right)  -1}-2^{\left(  m-2\right)  -1}=2^{\left(
m-2\right)  -1}\geq0$). Hence, $2^{\left(  m-2\right)  -1}\leq2^{\left(
m-1\right)  -1}$.

But the recursive definition of the Fibonacci sequence yields $f_{m}%
=f_{m-1}+f_{m-2}$ (since $m\geq2$). Hence,%
\[
f_{m}=\underbrace{f_{m-1}}_{\leq2^{\left(  m-1\right)  -1}}%
+\underbrace{f_{m-2}}_{\leq2^{\left(  m-2\right)  -1}\leq2^{\left(
m-1\right)  -1}}\leq2^{\left(  m-1\right)  -1}+2^{\left(  m-1\right)
-1}=2\cdot2^{\left(  m-1\right)  -1}=2^{m-1}.
\]
In other words, the statement $\mathcal{A}\left(  m\right)  $ holds (since the
statement $\mathcal{A}\left(  m\right)  $ is defined to be $\left(  f_{m}%
\leq2^{m-1}\right)  $).

Now, forget that we fixed $m$. We thus have shown that if $m\in\mathbb{Z}%
_{\geq0}$ is such that (\ref{pf.prop.ind.sip-fib.2nd.IH}) holds, then
$\mathcal{A}\left(  m\right)  $ holds. This completes the induction step.
Hence, by strong induction, we conclude that $\mathcal{A}\left(  n\right)  $
holds for each $n\in\mathbb{Z}_{\geq0}$. This completes our proof of
Proposition \ref{prop.ind.sip-fib}.
\end{proof}

The proof that we just showed still has a lot of \textquotedblleft
boilerplate\textquotedblright\ text that conveys no information. For example,
we have again explicitly defined the statement $\mathcal{A}\left(  n\right)
$, which is unnecessary: This statement is exactly what one would expect
(namely, the claim that we are proving, without the \textquotedblleft for each
$n\in\mathbb{N}$\textquotedblright\ part). Thus, in our case, this statement
is simply (\ref{eq.prop.ind.sip-fib.claim}). Furthermore, we can remove the
two sentences

\begin{quote}
\textquotedblleft Now, forget that we fixed $m$. We thus have shown that if
$m\in\mathbb{Z}_{\geq0}$ is such that (\ref{pf.prop.ind.sip-fib.2nd.IH})
holds, then $\mathcal{A}\left(  m\right)  $ holds.\textquotedblright.
\end{quote}

\noindent In fact, these sentences merely say that we have completed the
induction step; but this is clear anyway when we say that the induction step
is completed.

We said that we are proving our statement \textquotedblleft by strong
induction on $n$ starting at $0$\textquotedblright. Again, we can omit the
words \textquotedblleft starting at $0$\textquotedblright\ here, since this is
the only option (because our statement is about all $n\in\mathbb{Z}_{\geq0}$).

Finally, we can remove the words \textquotedblleft\textit{Induction
step:}\textquotedblright, because a proof by strong induction (unlike a proof
by standard induction) does not have an induction base (so the induction step
is all that it consists of).

Thus, our above proof can be shortened to the following:

\begin{proof}
[Proof of Proposition \ref{prop.ind.sip-fib} (third version).]We shall prove
(\ref{eq.prop.ind.sip-fib.claim}) by strong induction on $n$:

Let $m\in\mathbb{Z}_{\geq0}$. Assume that (\ref{eq.prop.ind.sip-fib.claim})
holds for every $n\in\mathbb{Z}_{\geq0}$ satisfying $n<m$. We must then show
that (\ref{eq.prop.ind.sip-fib.claim}) holds for $n=m$. In other words, we
must show that $f_{m}\leq2^{m-1}$ holds.

\begin{vershort}
This is true if $m=0$ (since $f_{0}=0\leq2^{0-1}$) and also true if $m=1$
(since $f_{1}=1=2^{1-1}$ and thus $f_{1}\leq2^{1-1}$). In other words, this is
true if $m\in\left\{  0,1\right\}  $. Thus, for the rest of the induction
step, we WLOG assume that we don't have $m\in\left\{  0,1\right\}  $. Hence,
$m\notin\left\{  0,1\right\}  $, so that $m\in\mathbb{N}\setminus\left\{
0,1\right\}  =\left\{  2,3,4,\ldots\right\}  $. Hence, $m\geq2$.
\end{vershort}

\begin{verlong}
This is true if $m=0$ (since $f_{0}=0\leq2^{0-1}$). Hence, for the rest of the
induction step, we WLOG assume that we don't have $m=0$. Hence, $m\neq0$, so
that $m\geq1$ (since $m\in\mathbb{Z}_{\geq0}$).

We must prove that $f_{m}\leq2^{m-1}$. This is true if $m=1$ (because
$f_{1}=1=2^{1-1}$ and thus $f_{1}\leq2^{1-1}$). Thus, for the rest of the
induction step, we WLOG assume that we don't have $m=1$. Hence, $m\neq1$.
Combining this with $m\geq1$, we obtain $m>1$, so that $m\geq2$ (since $m$ is
an integer).
\end{verlong}

From $m\geq2$, we conclude that $m-1\geq2-1=1\geq0$ and $m-2\geq2-2=0$. Thus,
both $m-1$ and $m-2$ belong to $\mathbb{N}$; therefore, $f_{m-1}$ and
$f_{m-2}$ are well-defined.

We have $m-1\in\mathbb{N}=\mathbb{Z}_{\geq0}$ and $m-1<m$. Hence,
(\ref{eq.prop.ind.sip-fib.claim}) (applied to $n=m-1$) yields that
$f_{m-1}\leq2^{\left(  m-1\right)  -1}$ (since we have assumed that
(\ref{eq.prop.ind.sip-fib.claim}) holds for every $n\in\mathbb{Z}_{\geq0}$
satisfying $n<m$).

We have $m-2\in\mathbb{N}=\mathbb{Z}_{\geq0}$ and $m-2<m$. Hence,
(\ref{eq.prop.ind.sip-fib.claim}) (applied to $n=m-2$) yields that
$f_{m-2}\leq2^{\left(  m-2\right)  -1}$ (since we have assumed that
(\ref{eq.prop.ind.sip-fib.claim}) holds for every $n\in\mathbb{Z}_{\geq0}$
satisfying $n<m$).

We have $\left(  m-1\right)  -1=m-2$ and thus $2^{\left(  m-1\right)
-1}=2^{m-2}=2\cdot2^{\left(  m-2\right)  -1}\geq2^{\left(  m-2\right)  -1}$
(since $2\cdot2^{\left(  m-2\right)  -1}-2^{\left(  m-2\right)  -1}=2^{\left(
m-2\right)  -1}\geq0$). Hence, $2^{\left(  m-2\right)  -1}\leq2^{\left(
m-1\right)  -1}$.

But the recursive definition of the Fibonacci sequence yields $f_{m}%
=f_{m-1}+f_{m-2}$ (since $m\geq2$). Hence,%
\[
f_{m}=\underbrace{f_{m-1}}_{\leq2^{\left(  m-1\right)  -1}}%
+\underbrace{f_{m-2}}_{\leq2^{\left(  m-2\right)  -1}\leq2^{\left(
m-1\right)  -1}}\leq2^{\left(  m-1\right)  -1}+2^{\left(  m-1\right)
-1}=2\cdot2^{\left(  m-1\right)  -1}=2^{m-1}.
\]
In other words, (\ref{eq.prop.ind.sip-fib.claim}) holds for $n=m$. This
completes the induction step. Hence, by strong induction, we conclude that
(\ref{eq.prop.ind.sip-fib.claim}) holds for each $n\in\mathbb{Z}_{\geq0}$. In
other words, (\ref{eq.prop.ind.sip-fib.claim}) holds for each $n\in\mathbb{N}$
(since $\mathbb{Z}_{\geq0}=\mathbb{N}$). This completes our proof of
Proposition \ref{prop.ind.sip-fib}.
\end{proof}

\subsection{Two unexpected integralities}

\subsubsection{The first integrality}

We shall illustrate strong induction on two further examples, which both have
the form of an \textquotedblleft unexpected integrality\textquotedblright: A
sequence of rational numbers is defined recursively by an equation that
involves fractions, but it turns out that all the entries of the sequence are
nevertheless integers. There is by now a whole genre of such results (see
\cite[Chapter 1]{Gale98} for an introduction\footnote{See also \cite{FomZel02}
for a seminal research paper at a more advanced level.}), and many of them are
connected with recent research in the theory of cluster algebras (see
\cite{Lampe} for an introduction).

The first of these examples is the following result:

\begin{proposition}
\label{prop.ind.LP1}Define a sequence $\left(  t_{0},t_{1},t_{2}%
,\ldots\right)  $ of positive rational numbers recursively by setting%
\begin{align*}
t_{0}  &  =1,\ \ \ \ \ \ \ \ \ \ t_{1}=1,\ \ \ \ \ \ \ \ \ \ t_{2}%
=1,\ \ \ \ \ \ \ \ \ \ \text{and}\\
t_{n}  &  =\dfrac{1+t_{n-1}t_{n-2}}{t_{n-3}}\ \ \ \ \ \ \ \ \ \ \text{for each
}n\geq3.
\end{align*}
(Thus,
\begin{align*}
t_{3}  &  =\dfrac{1+t_{2}t_{1}}{t_{0}}=\dfrac{1+1\cdot1}{1}=2;\\
t_{4}  &  =\dfrac{1+t_{3}t_{2}}{t_{1}}=\dfrac{1+2\cdot1}{1}=3;\\
t_{5}  &  =\dfrac{1+t_{4}t_{3}}{t_{2}}=\dfrac{1+3\cdot2}{1}=7;\\
t_{6}  &  =\dfrac{1+t_{5}t_{4}}{t_{3}}=\dfrac{1+7\cdot3}{2}=11,
\end{align*}
and so on.) Then:

\textbf{(a)} We have $t_{n+2}=4t_{n}-t_{n-2}$ for each $n\in\mathbb{Z}_{\geq
2}$.

\textbf{(b)} We have $t_{n}\in\mathbb{N}$ for each $n\in\mathbb{N}$.
\end{proposition}

Note that the sequence $\left(  t_{0},t_{1},t_{2},\ldots\right)  $ in
Proposition \ref{prop.ind.LP1} is clearly well-defined, because the expression
$\dfrac{1+t_{n-1}t_{n-2}}{t_{n-3}}$ always yields a well-defined positive
rational number when $t_{n-1},t_{n-2},t_{n-3}$ are positive rational numbers.
(In particular, the denominator $t_{n-3}$ of this fraction is $\neq0$ because
it is positive.) In contrast, if we had set $t_{n}=\dfrac{1-t_{n-1}t_{n-2}%
}{t_{n-3}}$ instead of $t_{n}=\dfrac{1+t_{n-1}t_{n-2}}{t_{n-3}}$, then the
sequence would \textbf{not} be well-defined (because then, we would get
$t_{3}=\dfrac{1-1\cdot1}{1}=0$ and $t_{6}=\dfrac{1-t_{5}t_{4}}{t_{3}}%
=\dfrac{1-t_{5}t_{4}}{0}$, which is undefined).

\begin{remark}
The sequence $\left(  t_{0},t_{1},t_{2},\ldots\right)  $ defined in
Proposition \ref{prop.ind.LP1} is the sequence
\href{https://oeis.org/A005246}{A005246 in the OEIS (Online Encyclopedia of
Integer Sequences)}. Its first entries are%
\begin{align*}
t_{0}  &  =1,\qquad t_{1}=1,\qquad t_{2}=1,\qquad t_{3}=2,\qquad t_{4}=3,\\
t_{5}  &  =7,\qquad t_{6}=11,\qquad t_{7}=26,\qquad t_{8}=41,\qquad t_{9}=97.
\end{align*}

Proposition \ref{prop.ind.LP1} \textbf{(b)} is an instance of the
\textit{Laurent phenomenon} (see, e.g., \cite[Example 3.2]{FomZel02}).
\end{remark}

Part \textbf{(a)} of Proposition \ref{prop.ind.LP1} is proven by a (regular)
induction; it is part \textbf{(b)} where strong induction comes handy:

\begin{proof}
[Proof of Proposition \ref{prop.ind.LP1}.]First, we notice that the recursive
definition of the sequence $\left(  t_{0},t_{1},t_{2},\ldots\right)  $ yields%
\begin{align*}
t_{3}  &  =\dfrac{1+t_{3-1}t_{3-2}}{t_{3-3}}=\dfrac{1+t_{2}t_{1}}{t_{0}%
}=\dfrac{1+1\cdot1}{1}\ \ \ \ \ \ \ \ \ \ \left(  \text{since }t_{0}=1\text{
and }t_{1}=1\text{ and }t_{2}=1\right) \\
&  =2.
\end{align*}
Furthermore, the recursive definition of the sequence $\left(  t_{0}%
,t_{1},t_{2},\ldots\right)  $ yields%
\begin{align*}
t_{4}  &  =\dfrac{1+t_{4-1}t_{4-2}}{t_{4-3}}=\dfrac{1+t_{3}t_{2}}{t_{1}%
}=\dfrac{1+2\cdot1}{1}\ \ \ \ \ \ \ \ \ \ \left(  \text{since }t_{1}=1\text{
and }t_{2}=1\text{ and }t_{3}=2\right) \\
&  =3.
\end{align*}
Thus, $t_{2+2}=t_{4}=3$. Comparing this with $4\underbrace{t_{2}}%
_{=1}-\underbrace{t_{2-2}}_{=t_{0}=1}=4\cdot1-1=3$, we obtain $t_{2+2}%
=4t_{2}-t_{2-2}$.

\textbf{(a)} We shall prove Proposition \ref{prop.ind.LP1} \textbf{(a)} by
induction on $n$ starting at $2$:

\textit{Induction base:} We have already shown that $t_{2+2}=4t_{2}-t_{2-2}$.
In other words, Proposition \ref{prop.ind.LP1} \textbf{(a)} holds for $n=2$.
This completes the induction base.

\textit{Induction step:} Let $m\in\mathbb{Z}_{\geq2}$. Assume that Proposition
\ref{prop.ind.LP1} \textbf{(a)} holds for $n=m$. We must prove that
Proposition \ref{prop.ind.LP1} \textbf{(a)} holds for $n=m+1$.

We have assumed that Proposition \ref{prop.ind.LP1} \textbf{(a)} holds for
$n=m$. In other words, we have $t_{m+2}=4t_{m}-t_{m-2}$.

We have $m\in\mathbb{Z}_{\geq2}$. Thus, $m$ is an integer that is $\geq2$.
Hence, $m\geq2$ and thus $m+1\geq2+1=3$. Thus, the recursive definition of the
sequence $\left(  t_{0},t_{1},t_{2},\ldots\right)  $ yields%
\[
t_{m+1}=\dfrac{1+t_{\left(  m+1\right)  -1}t_{\left(  m+1\right)  -2}%
}{t_{\left(  m+1\right)  -3}}=\dfrac{1+t_{m}t_{m-1}}{t_{m-2}}.
\]
Multiplying this equality by $t_{m-2}$, we obtain $t_{m-2}t_{m+1}%
=1+t_{m}t_{m-1}$. In other words,
\begin{equation}
t_{m-2}t_{m+1}-1=t_{m}t_{m-1}. \label{pf.prop.ind.LP1.a.3}%
\end{equation}
Hence,%
\begin{align}
1+\underbrace{t_{m+2}}_{=4t_{m}-t_{m-2}}t_{m+1}  &  =1+\left(  4t_{m}%
-t_{m-2}\right)  t_{m+1}=4t_{m}t_{m+1}-\underbrace{\left(  t_{m-2}%
t_{m+1}-1\right)  }_{\substack{=t_{m}t_{m-1}\\\text{(by
(\ref{pf.prop.ind.LP1.a.3}))}}}\nonumber\\
&  =4t_{m}t_{m+1}-t_{m}t_{m-1}=t_{m}\left(  4t_{m+1}-t_{m-1}\right)  .
\label{pf.prop.ind.LP1.a.4}%
\end{align}

Also, $m+3\geq3$. Thus, the recursive definition of the sequence $\left(
t_{0},t_{1},t_{2},\ldots\right)  $ yields%
\begin{align*}
t_{m+3}  &  =\dfrac{1+t_{\left(  m+3\right)  -1}t_{\left(  m+3\right)  -2}%
}{t_{\left(  m+3\right)  -3}}=\dfrac{1+t_{m+2}t_{m+1}}{t_{m}}=\dfrac{1}{t_{m}%
}\underbrace{\left(  1+t_{m+2}t_{m+1}\right)  }_{\substack{=t_{m}\left(
4t_{m+1}-t_{m-1}\right)  \\\text{(by (\ref{pf.prop.ind.LP1.a.4}))}}}\\
&  =\dfrac{1}{t_{m}}t_{m}\left(  4t_{m+1}-t_{m-1}\right)  =4t_{m+1}-t_{m-1}.
\end{align*}
In view of $m+3=\left(  m+1\right)  +2$ and $m-1=\left(  m+1\right)  -2$, this
rewrites as $t_{\left(  m+1\right)  +2}=4t_{m+1}-t_{\left(  m+1\right)  -2}$.
In other words, Proposition \ref{prop.ind.LP1} \textbf{(a)} holds for $n=m+1$.
This completes the induction step. Hence, Proposition \ref{prop.ind.LP1}
\textbf{(a)} is proven by induction.

\textbf{(b)} We shall prove Proposition \ref{prop.ind.LP1} \textbf{(b)} by
strong induction on $n$ starting at $0$:

\textit{Induction step:} Let $m\in\mathbb{N}$.\ \ \ \ \footnote{In order to
match the notations used in Theorem \ref{thm.ind.SIP}, we should be saying
\textquotedblleft Let $m\in\mathbb{Z}_{\geq0}$\textquotedblright\ here, rather
than \textquotedblleft Let $m\in\mathbb{N}$\textquotedblright. But of course,
this amounts to the same thing, since $\mathbb{N}=\mathbb{Z}_{\geq0}$.} Assume
that Proposition \ref{prop.ind.LP1} \textbf{(b)} holds for every
$n\in\mathbb{N}$ satisfying $n<m$. We must now show that Proposition
\ref{prop.ind.LP1} \textbf{(b)} holds for $n=m$.

We have assumed that Proposition \ref{prop.ind.LP1} \textbf{(b)} holds for
every $n\in\mathbb{N}$ satisfying $n<m$. In other words, we have
\begin{equation}
t_{n}\in\mathbb{N}\text{ for every }n\in\mathbb{N}\text{ satisfying }n<m.
\label{pf.prop.ind.LP1.b.IH}%
\end{equation}

We must now show that Proposition \ref{prop.ind.LP1} \textbf{(b)} holds for
$n=m$. In other words, we must show that $t_{m}\in\mathbb{N}$.

Recall that $\left(  t_{0},t_{1},t_{2},\ldots\right)  $ is a sequence of
positive rational numbers. Thus, $t_{m}$ is a positive rational number.

We are in one of the following five cases:

\textit{Case 1:} We have $m=0$.

\textit{Case 2:} We have $m=1$.

\textit{Case 3:} We have $m=2$.

\textit{Case 4:} We have $m=3$.

\textit{Case 5:} We have $m>3$.

Let us first consider Case 1. In this case, we have $m=0$. Thus, $t_{m}%
=t_{0}=1\in\mathbb{N}$. Hence, $t_{m}\in\mathbb{N}$ is proven in Case 1.

Similarly, we can prove $t_{m}\in\mathbb{N}$ in Case 2 (using $t_{1}=1$) and
in Case 3 (using $t_{2}=1$) and in Case 4 (using $t_{3}=2$). It thus remains
to prove $t_{m}\in\mathbb{N}$ in Case 5.

So let us consider Case 5. In this case, we have $m>3$. Thus, $m\geq4$ (since
$m$ is an integer), so that $m-2\geq4-2=2$. Thus, $m-2$ is an integer that is
$\geq2$. In other words, $m-2\in\mathbb{Z}_{\geq2}$. Hence, Proposition
\ref{prop.ind.LP1} \textbf{(a)} (applied to $n=m-2$) yields $t_{\left(
m-2\right)  +2}=4t_{m-2}-t_{\left(  m-2\right)  -2}$. In view of $\left(
m-2\right)  +2=m$ and $\left(  m-2\right)  -2=m-4$, this rewrites as
$t_{m}=4t_{m-2}-t_{m-4}$.

\begin{vershort}
But $m\geq4$, so that $m-4\in\mathbb{N}$, and $m-4<m$. Hence,
(\ref{pf.prop.ind.LP1.b.IH}) (applied to $n=m-4$) yields $t_{m-4}\in
\mathbb{N}\subseteq\mathbb{Z}$. Similarly, $t_{m-2}\in\mathbb{Z}$.
\end{vershort}

\begin{verlong}
But $m-2\in\mathbb{N}$ (since $m\geq4\geq2$) and $m-2<m$. Hence,
(\ref{pf.prop.ind.LP1.b.IH}) (applied to $n=m-2$) yields $t_{m-2}\in
\mathbb{N}\subseteq\mathbb{Z}$.

Also, $m-4\in\mathbb{N}$ (since $m\geq4$) and $m-4<m$. Hence,
(\ref{pf.prop.ind.LP1.b.IH}) (applied to $n=m-4$) yields $t_{m-4}\in
\mathbb{N}\subseteq\mathbb{Z}$.
\end{verlong}

So we know that $t_{m-2}$ and $t_{m-4}$ are both integers (since $t_{m-2}%
\in\mathbb{Z}$ and $t_{m-4}\in\mathbb{Z}$). Hence, $4t_{m-2}-t_{m-4}$ is an
integer as well. In other words, $t_{m}$ is an integer (because $t_{m}%
=4t_{m-2}-t_{m-4}$). Since $t_{m}$ is positive, we thus conclude that $t_{m}$
is a positive integer. Hence, $t_{m}\in\mathbb{N}$. This shows that $t_{m}%
\in\mathbb{N}$ in Case 5.

We now have proven $t_{m}\in\mathbb{N}$ in each of the five Cases 1, 2, 3, 4
and 5. Since these five Cases cover all possibilities, we thus conclude that
$t_{m}\in\mathbb{N}$ always holds. In other words, Proposition
\ref{prop.ind.LP1} \textbf{(b)} holds for $n=m$. This completes the induction
step. Thus, Proposition \ref{prop.ind.LP1} \textbf{(b)} is proven by strong induction.
\end{proof}

\subsubsection{The second integrality}

Our next example of an \textquotedblleft unexpected
integrality\textquotedblright\ is the following fact:

\begin{proposition}
\label{prop.ind.LP2}Fix a positive integer $r$. Define a sequence $\left(
b_{0},b_{1},b_{2},\ldots\right)  $ of positive rational numbers recursively by
setting%
\begin{align*}
b_{0}  &  =1,\ \ \ \ \ \ \ \ \ \ b_{1}=1,\ \ \ \ \ \ \ \ \ \ \text{and}\\
b_{n}  &  =\dfrac{b_{n-1}^{r}+1}{b_{n-2}}\ \ \ \ \ \ \ \ \ \ \text{for each
}n\geq2.
\end{align*}
(Thus,%
\begin{align*}
b_{2}  &  =\dfrac{b_{1}^{r}+1}{b_{0}}=\dfrac{1^{r}+1}{1}=2;\\
b_{3}  &  =\dfrac{b_{2}^{r}+1}{b_{1}}=\dfrac{2^{r}+1}{1}=2^{r}+1;\\
b_{4}  &  =\dfrac{b_{3}^{r}+1}{b_{2}}=\dfrac{\left(  2^{r}+1\right)  ^{r}%
+1}{2},
\end{align*}
and so on.) Then:

\textbf{(a)} We have $b_{n}\in\mathbb{N}$ for each $n\in\mathbb{N}$.

\textbf{(b)} If $r\geq2$, then $b_{n}\mid b_{n-2}+b_{n+2}$ for each
$n\in\mathbb{Z}_{\geq2}$.
\end{proposition}

\begin{remark}
If $r=1$, then the sequence $\left(  b_{0},b_{1},b_{2},\ldots\right)  $
defined in Proposition \ref{prop.ind.LP2} is
\[
\left(  1,1,2,3,2,1,1,2,3,2,1,1,2,3,2,\ldots\right)
\]
(this is a periodic sequence, which consists of the five terms $1,1,2,3,2$
repeated over and over); this can easily be proven by strong induction.
Despite its simplicity, this sequence is the sequence
\href{https://oeis.org/A076839}{A076839 in the OEIS}.

If $r=2$, then the sequence $\left(  b_{0},b_{1},b_{2},\ldots\right)  $
defined in Proposition \ref{prop.ind.LP2} is%
\[
\left(  1,f_{1},f_{3},f_{5},f_{7},\ldots\right)  =\left(
1,1,2,5,13,34,89,233,610,1597,\ldots\right)
\]
consisting of all Fibonacci numbers at odd positions (i.e., Fibonacci numbers
of the form $f_{2n-1}$ for $n\in\mathbb{N}$) with an extra $1$ at the front.
This, again, can be proven by induction. Also, this sequence satisfies the
recurrence relation $b_{n}=3b_{n-1}-b_{n-2}$ for all $n\geq2$. This is the
sequence \href{https://oeis.org/A001519}{A001519 in the OEIS}.

If $r=3$, then the sequence $\left(  b_{0},b_{1},b_{2},\ldots\right)  $
defined in Proposition \ref{prop.ind.LP2} is%
\[
\left(  1,1,2,9,365,5403014,432130991537958813,\ldots\right)  ;
\]
its entries grow so fast that the next entry would need a separate line. This
is the sequence \href{https://oeis.org/A003818}{A003818 in the OEIS}. Unlike
the cases of $r=1$ and $r=2$, not much can be said about this sequence, other
than what has been said in Proposition \ref{prop.ind.LP2}.

Proposition \ref{prop.ind.LP2} \textbf{(a)} is an instance of the
\textit{Laurent phenomenon for cluster algebras} (see, e.g., \cite[Example
2.5]{FomZel01}; also, see \cite{Marsh13} and \cite{FoWiZe16} for expositions).
See also \cite{MusPro07} for a study of the specific recurrence equation from
Proposition \ref{prop.ind.LP2} (actually, a slightly more general equation).
\end{remark}

Before we prove Proposition \ref{prop.ind.LP2}, let us state an auxiliary fact:

\begin{lemma}
\label{lem.ind.LP2.Hx}Let $r\in\mathbb{N}$. For every nonzero $x\in\mathbb{Q}%
$, we set $H\left(  x\right)  =\dfrac{\left(  x+1\right)  ^{r}-1}{x}$.

Then, $H\left(  x\right)  \in\mathbb{Z}$ whenever $x$ is a nonzero integer.
\end{lemma}

\begin{proof}
[Proof of Lemma \ref{lem.ind.LP2.Hx}.]Let $x$ be a nonzero integer. Then,
$x\mid\left(  x+1\right)  -1$ (because $\left(  x+1\right)  -1=x$). In other
words, $x+1\equiv1\operatorname{mod}x$ (by the definition of \textquotedblleft
congruent\textquotedblright). Hence, Proposition \ref{prop.mod.pow} (applied
to $a=x+1$, $b=1$, $n=x$ and $k=r$) shows that $\left(  x+1\right)  ^{r}%
\equiv1^{r}=1\operatorname{mod}x$. In other words, $\left(  x+1\right)
^{r}-1$ is divisible by $x$. In other words, $\dfrac{\left(  x+1\right)
^{r}-1}{x}$ is an integer. In other words, $\dfrac{\left(  x+1\right)  ^{r}%
-1}{x}\in\mathbb{Z}$. Thus, $H\left(  x\right)  =\dfrac{\left(  x+1\right)
^{r}-1}{x}\in\mathbb{Z}$. This proves Lemma \ref{lem.ind.LP2.Hx}.
\end{proof}

\begin{proof}
[Proof of Proposition \ref{prop.ind.LP2}.]First, we notice that the recursive
definition of the sequence $\left(  b_{0},b_{1},b_{2},\ldots\right)  $ yields%
\begin{align*}
b_{2}  &  =\dfrac{b_{2-1}^{r}+1}{b_{2-2}}=\dfrac{b_{1}^{r}+1}{b_{0}}%
=\dfrac{1^{r}+1}{1}\ \ \ \ \ \ \ \ \ \ \left(  \text{since }b_{0}=1\text{ and
}b_{1}=1\right) \\
&  =\dfrac{1+1}{1}\ \ \ \ \ \ \ \ \ \ \left(  \text{since }1^{r}=1\right) \\
&  =2.
\end{align*}
Furthermore, the recursive definition of the sequence $\left(  b_{0}%
,b_{1},b_{2},\ldots\right)  $ yields%
\begin{align*}
b_{3}  &  =\dfrac{b_{3-1}^{r}+1}{b_{3-2}}=\dfrac{b_{2}^{r}+1}{b_{1}}%
=\dfrac{2^{r}+1}{1}\ \ \ \ \ \ \ \ \ \ \left(  \text{since }b_{1}=1\text{ and
}b_{2}=2\right) \\
&  =2^{r}+1.
\end{align*}

For every nonzero $x\in\mathbb{Q}$, we set $H\left(  x\right)  =\dfrac{\left(
x+1\right)  ^{r}-1}{x}$.

For every integer $m\geq1$, we have%
\begin{equation}
b_{m}^{r}+1=b_{m+1}b_{m-1}. \label{pf.prop.ind.LP2.o0}%
\end{equation}

[\textit{Proof of (\ref{pf.prop.ind.LP2.o0}):} Let $m\geq1$ be an integer.
From $m\geq1$, we obtain $m+1\geq1+1=2$. Hence, the recursive definition of
the sequence $\left(  b_{0},b_{1},b_{2},\ldots\right)  $ yields
\[
b_{m+1}=\dfrac{b_{\left(  m+1\right)  -1}^{r}+1}{b_{\left(  m+1\right)  -2}%
}=\dfrac{b_{m}^{r}+1}{b_{m-1}}.
\]
Multiplying both sides of this equality by $b_{m-1}$, we obtain $b_{m+1}%
b_{m-1}=b_{m}^{r}+1$. This proves (\ref{pf.prop.ind.LP2.o0}).]

Let us first prove the following observation:

\begin{statement}
\textit{Observation 1:} Each integer $n\geq2$ satisfies $b_{n+2}%
=b_{n-2}b_{n+1}^{r}-b_{n}^{r-1}H\left(  b_{n}^{r}\right)  $.
\end{statement}

[\textit{Proof of Observation 1:} Let $n\geq2$ be an integer. Thus,
$n\geq2\geq1$. Thus, (\ref{pf.prop.ind.LP2.o0}) (applied to $m=n$) yields%
\begin{equation}
b_{n}^{r}+1=b_{n+1}b_{n-1}. \label{pf.prop.ind.LP2.o1.pf.1}%
\end{equation}

On the other hand, $n+1\geq n\geq2\geq1$. Hence, (\ref{pf.prop.ind.LP2.o0})
(applied to $m=n+1$) yields%
\[
b_{n+1}^{r}+1=b_{\left(  n+1\right)  +1}b_{\left(  n+1\right)  -1}%
=b_{n+2}b_{n}.
\]
Hence,%
\begin{equation}
b_{n+1}^{r}=b_{n+2}b_{n}-1. \label{pf.prop.ind.LP2.o1.pf.2}%
\end{equation}

Also, $n-1\geq1$ (since $n\geq2=1+1$). Hence, (\ref{pf.prop.ind.LP2.o0})
(applied to $m=n-1$) yields%
\[
b_{n-1}^{r}+1=b_{\left(  n-1\right)  +1}b_{\left(  n-1\right)  -1}%
=b_{n}b_{n-2}.
\]
Hence,%
\begin{equation}
b_{n-1}^{r}=b_{n}b_{n-2}-1. \label{pf.prop.ind.LP2.o1.pf.3}%
\end{equation}

But $b_{n}$ is a positive rational number (since $\left(  b_{0},b_{1}%
,b_{2},\ldots\right)  $ is a sequence of positive rational numbers). Thus,
$b_{n}^{r}$ is also a positive rational number. Hence, $b_{n}^{r}\in
\mathbb{Q}$ is nonzero. The definition of $H\left(  b_{n}^{r}\right)  $ yields
$H\left(  b_{n}^{r}\right)  =\dfrac{\left(  b_{n}^{r}+1\right)  ^{r}-1}%
{b_{n}^{r}}$; therefore,
\begin{align*}
b_{n}^{r-1}\underbrace{H\left(  b_{n}^{r}\right)  }_{=\dfrac{\left(  b_{n}%
^{r}+1\right)  ^{r}-1}{b_{n}^{r}}}  &  =b_{n}^{r-1}\cdot\dfrac{\left(
b_{n}^{r}+1\right)  ^{r}-1}{b_{n}^{r}}=\underbrace{\dfrac{b_{n}^{r-1}}%
{b_{n}^{r}}}_{=\dfrac{1}{b_{n}}}\cdot\left(  \left(  \underbrace{b_{n}^{r}%
+1}_{\substack{=b_{n+1}b_{n-1}\\\text{(by (\ref{pf.prop.ind.LP2.o1.pf.1}))}%
}}\right)  ^{r}-1\right) \\
&  =\dfrac{1}{b_{n}}\cdot\left(  \underbrace{\left(  b_{n+1}b_{n-1}\right)
^{r}}_{=b_{n+1}^{r}b_{n-1}^{r}}-1\right)  =\dfrac{1}{b_{n}}\cdot\left(
\underbrace{b_{n+1}^{r}}_{\substack{=b_{n+2}b_{n}-1\\\text{(by
(\ref{pf.prop.ind.LP2.o1.pf.2}))}}}\ \ \underbrace{b_{n-1}^{r}}%
_{\substack{=b_{n}b_{n-2}-1\\\text{(by (\ref{pf.prop.ind.LP2.o1.pf.3}))}%
}}-1\right) \\
&  =\dfrac{1}{b_{n}}\cdot\underbrace{\left(  \left(  b_{n+2}b_{n}-1\right)
\left(  b_{n}b_{n-2}-1\right)  -1\right)  }_{=b_{n}\left(  b_{n}b_{n+2}%
b_{n-2}-b_{n+2}-b_{n-2}\right)  }\\
&  =\dfrac{1}{b_{n}}\cdot b_{n}\left(  b_{n}b_{n+2}b_{n-2}-b_{n+2}%
-b_{n-2}\right) \\
&  =b_{n}b_{n+2}b_{n-2}-b_{n+2}-b_{n-2}=b_{n-2}\underbrace{\left(
b_{n+2}b_{n}-1\right)  }_{\substack{=b_{n+1}^{r}\\\text{(by
(\ref{pf.prop.ind.LP2.o1.pf.2}))}}}-b_{n+2}\\
&  =b_{n-2}b_{n+1}^{r}-b_{n+2}.
\end{align*}
Solving this equation for $b_{n+2}$, we obtain $b_{n+2}=b_{n-2}b_{n+1}%
^{r}-b_{n}^{r-1}H\left(  b_{n}^{r}\right)  $. This proves Observation 1.]

\textbf{(a)} We shall prove Proposition \ref{prop.ind.LP2} \textbf{(a)} by
strong induction on $n$:

\textit{Induction step:} Let $m\in\mathbb{N}$. Assume that Proposition
\ref{prop.ind.LP2} \textbf{(a)} holds for every $n \in\mathbb{N}$ satisfying
$n<m$. We must now prove that Proposition \ref{prop.ind.LP2} \textbf{(a)}
holds for $n=m$.

We have assumed that Proposition \ref{prop.ind.LP2} \textbf{(a)} holds for
every $n \in\mathbb{N}$ satisfying $n<m$. In other words, we have%
\begin{equation}
b_{n}\in\mathbb{N}\text{ for every }n\in\mathbb{N}\text{ satisfying }n<m.
\label{pf.prop.ind.LP2.a.IH}%
\end{equation}

We must now show that Proposition \ref{prop.ind.LP2} \textbf{(a)} holds for
$n=m$. In other words, we must show that $b_{m}\in\mathbb{N}$.

Recall that $\left(  b_{0},b_{1},b_{2},\ldots\right)  $ is a sequence of
positive rational numbers. Thus, $b_{m}$ is a positive rational number.

We are in one of the following five cases:

\textit{Case 1:} We have $m=0$.

\textit{Case 2:} We have $m=1$.

\textit{Case 3:} We have $m=2$.

\textit{Case 4:} We have $m=3$.

\textit{Case 5:} We have $m>3$.

Let us first consider Case 1. In this case, we have $m=0$. Thus, $b_{m}%
=b_{0}=1\in\mathbb{N}$. Hence, $b_{m}\in\mathbb{N}$ is proven in Case 1.

Similarly, we can prove $b_{m}\in\mathbb{N}$ in Case 2 (using $b_{1}=1$) and
in Case 3 (using $b_{2}=2$) and in Case 4 (using $b_{3}=2^{r}+1$). It thus
remains to prove $b_{m}\in\mathbb{N}$ in Case 5.

So let us consider Case 5. In this case, we have $m>3$. Thus, $m\geq4$ (since
$m$ is an integer), so that $m-2\geq4-2=2$. Hence, Observation 1 (applied to
$n=m-2$) yields $b_{\left(  m-2\right)  +2}=b_{\left(  m-2\right)
-2}b_{\left(  m-2\right)  +1}^{r}-b_{m-2}^{r-1}H\left(  b_{m-2}^{r}\right)  $.
In view of $\left(  m-2\right)  +2=m$ and $\left(  m-2\right)  -2=m-4$ and
$\left(  m-2\right)  +1=m-1$, this rewrites as
\begin{equation}
b_{m}=b_{m-4}b_{m-1}^{r}-b_{m-2}^{r-1}H\left(  b_{m-2}^{r}\right)  .
\label{pf.prop.ind.LP2.a.bm=}%
\end{equation}

But $m-2\in\mathbb{N}$ (since $m\geq4\geq2$) and $m-2<m$. Hence,
(\ref{pf.prop.ind.LP2.a.IH}) (applied to $n=m-2$) yields $b_{m-2}\in
\mathbb{N}\subseteq\mathbb{Z}$. Also, $b_{m-2}$ is a positive rational number
(since $\left(  b_{0},b_{1},b_{2},\ldots\right)  $ is a sequence of positive
rational numbers) and thus a positive integer (since $b_{m-2}\in\mathbb{N}$),
hence a nonzero integer. Thus, $b_{m-2}^{r}$ is a nonzero integer as well.
Therefore, Lemma \ref{lem.ind.LP2.Hx} (applied to $x=b_{m-2}^{r}$) shows that
$H\left(  b_{m-2}^{r}\right)  \in\mathbb{Z}$. In other words, $H\left(
b_{m-2}^{r}\right)  $ is an integer. Also, $r-1\geq0$ (since $r\geq1$), and
thus $r-1\in\mathbb{N}$. Hence, $b_{m-2}^{r-1}$ is an integer (since $b_{m-2}$
is an integer).

Also, $m-4\in\mathbb{N}$ (since $m\geq4$) and $m-4<m$. Hence,
(\ref{pf.prop.ind.LP2.a.IH}) (applied to $n=m-4$) yields $b_{m-4}\in
\mathbb{N}\subseteq\mathbb{Z}$. In other words, $b_{m-4}$ is an integer.

\begin{vershort}
Similarly, $b_{m-1}$ is an integer. Thus, $b_{m-1}^{r}$ is an integer.
\end{vershort}

\begin{verlong}
Also, $m-1\in\mathbb{N}$ (since $m\geq4\geq1$) and $m-1<m$. Hence,
(\ref{pf.prop.ind.LP2.a.IH}) (applied to $n=m-1$) yields $b_{m-1}\in
\mathbb{N}\subseteq\mathbb{Z}$. Hence, $b_{m-1}^{r}\in\mathbb{Z}$ (since
$r\in\mathbb{N}$ (because $r\geq1$)). In other words, $b_{m-1}^{r}$ is an integer.
\end{verlong}

We now know that the four numbers $b_{m-4}$, $b_{m-1}^{r}$, $b_{m-2}^{r-1}$
and $H\left(  b_{m-2}^{r}\right)  $ are integers. Thus, the number
$b_{m-4}b_{m-1}^{r}-b_{m-2}^{r-1}H\left(  b_{m-2}^{r}\right)  $ also is an
integer (since it is obtained from these four numbers by multiplication and
subtraction). In view of (\ref{pf.prop.ind.LP2.a.bm=}), this rewrites as
follows: The number $b_{m}$ is an integer. Since $b_{m}$ is positive, we thus
conclude that $b_{m}$ is a positive integer. Hence, $b_{m}\in\mathbb{N}$. This
shows that $b_{m}\in\mathbb{N}$ in Case 5.

We now have proven $b_{m}\in\mathbb{N}$ in each of the five Cases 1, 2, 3, 4
and 5. Thus, $b_{m}\in\mathbb{N}$ always holds. In other words, Proposition
\ref{prop.ind.LP2} \textbf{(a)} holds for $n=m$. This completes the induction
step. Thus, Proposition \ref{prop.ind.LP2} \textbf{(a)} is proven by strong induction.

\textbf{(b)} Assume that $r\geq2$. We must prove that $b_{n}\mid
b_{n-2}+b_{n+2}$ for each $n\in\mathbb{Z}_{\geq2}$.

So let $n\in\mathbb{Z}_{\geq2}$. We must show that $b_{n}\mid b_{n-2}+b_{n+2}$.

\begin{vershort}
From $n\in\mathbb{Z}_{\geq2}$, we obtain $n\geq2$, so that $n-2\in\mathbb{N}$.
\end{vershort}

\begin{verlong}
We have $n\in\mathbb{Z}_{\geq2}$. Thus, $n$ is an integer that is $\geq2$.
Hence, $n\geq2$. Thus, $n-2\geq0$, so that $n-2\in\mathbb{N}$.
\end{verlong}

\begin{vershort}
Proposition \ref{prop.ind.LP2} \textbf{(a)} (applied to $n-2$ instead of $n$)
yields $b_{n-2}\in\mathbb{N}$. Similarly, $b_{n}\in\mathbb{N}$ and $b_{n+1}%
\in\mathbb{N}$ and $b_{n+2}\in\mathbb{N}$. Thus, all of $b_{n-2}$, $b_{n+1}$,
$b_{n}$ are $b_{n+2}$ are integers.
\end{vershort}

\begin{verlong}
Thus, Proposition \ref{prop.ind.LP2} \textbf{(a)} (applied to $n-2$ instead of
$n$) yields $b_{n-2}\in\mathbb{N}$. Also, $n\geq2$, so that $n\in\mathbb{N}$.
Hence, Proposition \ref{prop.ind.LP2} \textbf{(a)} yields $b_{n}\in\mathbb{N}%
$. Furthermore, $n+2\geq n\geq2$, so that $n+2\in\mathbb{N}$. Thus,
Proposition \ref{prop.ind.LP2} \textbf{(a)} (applied to $n+2$ instead of $n$)
yields $b_{n+2}\in\mathbb{N}$. Also, $n+1\geq n\geq2$, so that $n+1\in
\mathbb{N}$. Thus, Proposition \ref{prop.ind.LP2} \textbf{(a)} (applied to
$n+1$ instead of $n$) yields $b_{n+1}\in\mathbb{N}$. All of the numbers
$b_{n-2}$, $b_{n+1}$, $b_{n}$ are $b_{n+2}$ are integers (since $b_{n-2}%
\in\mathbb{N}\subseteq\mathbb{Z}$, $b_{n+1}\in\mathbb{N}\subseteq\mathbb{Z}$,
$b_{n}\in\mathbb{N}\subseteq\mathbb{Z}$ and $b_{n+2}\in\mathbb{N}%
\subseteq\mathbb{Z}$).
\end{verlong}

But $b_{n}$ is a positive rational number (since $\left(  b_{0},b_{1}%
,b_{2},\ldots\right)  $ is a sequence of positive rational numbers), and
therefore a positive integer (since $b_{n}$ is an integer). Hence, $b_{n}^{r}$
is a positive integer, and thus a nonzero integer. Therefore, Lemma
\ref{lem.ind.LP2.Hx} (applied to $x=b_{n}^{r}$) shows that $H\left(  b_{n}%
^{r}\right)  \in\mathbb{Z}$. In other words, $H\left(  b_{n}^{r}\right)  $ is
an integer.

We have $n+2\geq2$. Hence, the recursive definition of the sequence $\left(
b_{0},b_{1},b_{2},\ldots\right)  $ yields $b_{n+2}=\dfrac{b_{\left(
n+2\right)  -1}^{r}+1}{b_{\left(  n+2\right)  -2}}=\dfrac{b_{n+1}^{r}+1}%
{b_{n}}$. Multiplying this equality by $b_{n}$, we obtain%
\begin{equation}
b_{n}b_{n+2}=b_{n+1}^{r}+1. \label{pf.prop.ind.LP2.b.4}%
\end{equation}

\begin{vershort}
We have $r-2\in\mathbb{N}$ (since $r\geq2$) and $b_{n}\in\mathbb{N}$. Hence,
$b_{n}^{r-2}$ is an integer.
\end{vershort}

\begin{verlong}
We have $r\geq2$, so that $r-2\geq0$. Thus, $r-2\in\mathbb{N}$. Therefore,
$b_{n}^{r-2}$ is an integer (since $b_{n}$ is an integer).
\end{verlong}

Observation 1 yields $b_{n+2}=b_{n-2}b_{n+1}^{r}-b_{n}^{r-1}H\left(  b_{n}%
^{r}\right)  $. Thus,%
\begin{align}
&  \underbrace{b_{n+2}}_{=b_{n-2}b_{n+1}^{r}-b_{n}^{r-1}H\left(  b_{n}%
^{r}\right)  }+b_{n-2}\nonumber\\
&  =b_{n-2}b_{n+1}^{r}-b_{n}^{r-1}H\left(  b_{n}^{r}\right)  +b_{n-2}%
=\underbrace{b_{n-2}b_{n+1}^{r}+b_{n-2}}_{=b_{n-2}\left(  b_{n+1}%
^{r}+1\right)  }-b_{n}^{r-1}H\left(  b_{n}^{r}\right) \nonumber\\
&  =b_{n-2}\underbrace{\left(  b_{n+1}^{r}+1\right)  }_{\substack{=b_{n}%
b_{n+2}\\\text{(by (\ref{pf.prop.ind.LP2.b.4}))}}}-\underbrace{b_{n}^{r-1}%
}_{=b_{n}b_{n}^{r-2}}H\left(  b_{n}^{r}\right)  =b_{n-2}b_{n}b_{n+2}%
-b_{n}b_{n}^{r-2}H\left(  b_{n}^{r}\right) \nonumber\\
&  =b_{n}\left(  b_{n-2}b_{n+2}-b_{n}^{r-2}H\left(  b_{n}^{r}\right)  \right)
. \label{pf.prop.ind.LP2.b.9}%
\end{align}
But $b_{n-2}b_{n+2}-b_{n}^{r-2}H\left(  b_{n}^{r}\right)  $ is an integer
(because $b_{n-2}$, $b_{n+2}$, $b_{n}^{r-2}$ and $H\left(  b_{n}^{r}\right)  $
are integers). Denote this integer by $z$. Thus, $z=b_{n-2}b_{n+2}-b_{n}%
^{r-2}H\left(  b_{n}^{r}\right)  $. Since $b_{n}$ and $z$ are integers, we
have%
\begin{align*}
b_{n}  &  \mid b_{n}\underbrace{z}_{=b_{n-2}b_{n+2}-b_{n}^{r-2}H\left(
b_{n}^{r}\right)  }\\
&  =b_{n}\left(  b_{n-2}b_{n+2}-b_{n}^{r-2}H\left(  b_{n}^{r}\right)  \right)
=b_{n+2}+b_{n-2}\ \ \ \ \ \ \ \ \ \ \left(  \text{by
(\ref{pf.prop.ind.LP2.b.9})}\right) \\
&  =b_{n-2}+b_{n+2}.
\end{align*}
This proves Proposition \ref{prop.ind.LP2} \textbf{(b)}.
\end{proof}

For a (slightly) different proof of Proposition \ref{prop.ind.LP2}, see
\url{http://artofproblemsolving.com/community/c6h428645p3705719} .

\begin{remark}
You might wonder what happens if we replace \textquotedblleft$b_{n-1}^{r}%
+1$\textquotedblright\ by \textquotedblleft$b_{n-1}^{r}+q$\textquotedblright%
\ in Proposition \ref{prop.ind.LP2}, where $q$ is some fixed nonnegative
integer. The answer turns out to be somewhat disappointing in general: For
example, if we set $r=3$ and $q=2$, then our sequence $\left(  b_{0}%
,b_{1},b_{2},\ldots\right)  $ begins with%
\begin{align*}
b_{0} &  =1,\ \ \ \ \ \ \ \ \ \ b_{1}=1,\ \ \ \ \ \ \ \ \ \ b_{2}=\dfrac
{1^{3}+2}{1}=3,\\
b_{3} &  =\dfrac{3^{3}+2}{1}=29,\ \ \ \ \ \ \ \ \ \ b_{4}=\dfrac{29^{3}+2}%
{3}=\dfrac{24\,391}{3},
\end{align*}
at which point it becomes clear that $b_{n}\in\mathbb{N}$ no longer holds for
all $n\in\mathbb{N}$. The same happens for all $r>2$ as long as $q=3$. (More
precisely, if we take $r>2$ and $q=2$, then the first $n$ that violates
$b_{n}\in\mathbb{N}$ is $4$ or $5$ depending on whether $r$ is odd or even.
Proving this is a nice exercise!)

However, $b_{n}\in\mathbb{N}$ still holds for all $n\in\mathbb{N}$ when $r=2$.
This follows from Exercise \ref{exe.ind.LP2q} below.
\end{remark}

\begin{exercise}
\label{exe.ind.LP2q}Fix a nonnegative integer $q$. Define a sequence $\left(
b_{0},b_{1},b_{2},\ldots\right)  $ of positive rational numbers recursively by
setting%
\begin{align*}
b_{0}  &  =1,\ \ \ \ \ \ \ \ \ \ b_{1}=1,\ \ \ \ \ \ \ \ \ \ \text{and}\\
b_{n}  &  =\dfrac{b_{n-1}^{2}+q}{b_{n-2}}\ \ \ \ \ \ \ \ \ \ \text{for each
}n\geq2.
\end{align*}
(Thus,%
\begin{align*}
b_{2}  &  =\dfrac{b_{1}^{2}+q}{b_{0}}=\dfrac{1^{2}+q}{1}=q+1;\\
b_{3}  &  =\dfrac{b_{2}^{2}+q}{b_{1}}=\dfrac{\left(  q+1\right)  ^{2}+q}%
{1}=q^{2}+3q+1;\\
b_{4}  &  =\dfrac{b_{3}^{2}+q}{b_{2}}=\dfrac{\left(  q^{2}+3q+1\right)
^{2}+q}{q+1}=q^{3}+5q^{2}+6q+1;
\end{align*}
and so on.) Prove that:

\textbf{(a)} We have $b_{n}=\left(  q+2\right)  b_{n-1}-b_{n-2}$ for each
$n\in\mathbb{Z}_{\geq2}$.

\textbf{(b)} We have $b_{n}\in\mathbb{N}$ for each $n\in\mathbb{N}$.
\end{exercise}

\subsection{\label{sect.ind.bez}Strong induction on a derived quantity:
Bezout's theorem}

\subsubsection{Strong induction on a derived quantity}

In Section \ref{sect.ind.max}, we have seen how to use induction on a variable
that does not explicitly appear in the claim. In the current section, we shall
show the same for strong induction. This time, the fact that we shall be
proving is the following:

\begin{theorem}
\label{thm.ind.bezout}Let $a\in\mathbb{N}$ and $b\in\mathbb{N}$. Then, there
exist $g\in\mathbb{N}$, $x\in\mathbb{Z}$ and $y\in\mathbb{Z}$ such that
$g=ax+by$ and $g\mid a$ and $g\mid b$.
\end{theorem}

\begin{example}
\textbf{(a)} If $a=3$ and $b=5$, then Theorem \ref{thm.ind.bezout} says that
there exist $g\in\mathbb{N}$, $x\in\mathbb{Z}$ and $y\in\mathbb{Z}$ such that
$g=3x+5y$ and $g\mid3$ and $g\mid5$. And indeed, it is easy to find such $g$,
$x$ and $y$: for example, $g=1$, $x=-3$ and $y=2$ will do (since $1=3\left(
-3\right)  +5\cdot2$ and $1\mid3$ and $1\mid5$).

\textbf{(b)} If $a=4$ and $b=6$, then Theorem \ref{thm.ind.bezout} says that
there exist $g\in\mathbb{N}$, $x\in\mathbb{Z}$ and $y\in\mathbb{Z}$ such that
$g=4x+6y$ and $g\mid4$ and $g\mid6$. And indeed, it is easy to find such $g$,
$x$ and $y$: for example, $g=2$, $x=-1$ and $y=1$ will do (since $2=4\left(
-1\right)  +6\cdot1$ and $2\mid4$ and $2\mid6$).
\end{example}

Theorem \ref{thm.ind.bezout} is one form of \textit{Bezout's theorem for
integers}, and its real significance might not be clear at this point; it
becomes important when the greatest common divisor of two integers is studied.
For now, we observe that the $g$ in Theorem \ref{thm.ind.bezout} is obviously
a common divisor of $a$ and $b$ (that is, an integer that divides both $a$ and
$b$); but it is also divisible by every common divisor of $a$ and $b$ (because
of Proposition \ref{prop.div.ax+by}).

Let us now focus on the proof of Theorem \ref{thm.ind.bezout}. It is natural
to try proving it by induction (or perhaps strong induction) on $a$ or on $b$,
but neither option leads to success. It may feel like \textquotedblleft
induction on $a$ and on $b$ at the same time\textquotedblright\ could help,
and this is indeed a viable approach\footnote{Of course, it needs to be done
right: An induction proof always requires choosing \textbf{one} variable to do
induction on; but it is possible to nest an induction proof inside the
induction step (or inside the induction base) of a different induction proof.
For example, imagine that we are trying to prove that%
\[
ab=ba\ \ \ \ \ \ \ \ \ \ \text{for any }a\in\mathbb{N}\text{ and }%
b\in\mathbb{N}.
\]
We can prove this by induction on $a$. More precisely, for each $a\in
\mathbb{N}$, we let $\mathcal{A}\left(  a\right)  $ be the statement $\left(
ab=ba\text{ for all }b\in\mathbb{N}\right)  $. We then prove $\mathcal{A}%
\left(  a\right)  $ by induction on $a$. In the induction step, we fix
$m\in\mathbb{N}$, and we assume that $\mathcal{A}\left(  m\right)  $ holds; we
now need to prove that $\mathcal{A}\left(  m+1\right)  $ holds. In other
words, we need to prove that $\left(  m+1\right)  b=b\left(  m+1\right)  $ for
all $b\in\mathbb{N}$. We can now prove this statement by induction on $b$
(although there are easier options, of course). Thus, the induction proof of
this statement happens inside the induction step of another induction proof.
This nesting of induction proofs is legitimate (and even has a name: it is
called \textit{double induction}), but tends to be rather confusing (just
think about what the sentence \textquotedblleft The induction base is
complete\textquotedblright\ means: is it about the induction base of the first
induction proof, or that of the second?), and is best avoided when possible.}.
But there is a simpler and shorter method available: strong induction on
$a+b$. As in Section \ref{sect.ind.max}, the way to formalize such a strong
induction is by introducing auxiliary statements $\mathcal{A}\left(  n\right)
$, which say as much as \textquotedblleft Theorem \ref{thm.ind.bezout} holds
under the requirement that $a+b=n$\textquotedblright:

\begin{proof}
[Proof of Theorem \ref{thm.ind.bezout}.]First of all, let us forget that we
fixed $a$ and $b$. So we want to prove that if $a\in\mathbb{N}$ and
$b\in\mathbb{N}$, then there exist $g\in\mathbb{N}$, $x\in\mathbb{Z}$ and
$y\in\mathbb{Z}$ such that $g=ax+by$ and $g\mid a$ and $g\mid b$.

For each $n\in\mathbb{N}$, we let $\mathcal{A}\left(  n\right)  $ be the
statement%
\begin{equation}
\left(
\begin{array}
[c]{c}%
\text{if }a\in\mathbb{N}\text{ and }b\in\mathbb{N}\text{ satisfy }a+b=n\text{,
then there exist }g\in\mathbb{N}\text{, }x\in\mathbb{Z}\\
\text{and }y\in\mathbb{Z}\text{ such that }g=ax+by\text{ and }g\mid a\text{
and }g\mid b
\end{array}
\right)  . \label{pf.thm.ind.bezout.An}%
\end{equation}

We claim that $\mathcal{A}\left(  n\right)  $ holds for all $n\in\mathbb{N}$.

Indeed, let us prove this by strong induction on $n$ starting at $0$:

\textit{Induction step:} Let $m\in\mathbb{N}$. Assume that%
\begin{equation}
\left(  \mathcal{A}\left(  n\right)  \text{ holds for every }n\in
\mathbb{N}\text{ satisfying }n<m\right)  . \label{pf.thm.ind.bezout.IH}%
\end{equation}
We must then show that $\mathcal{A}\left(  m\right)  $ holds.

To do this, we shall prove the following claim:

\begin{statement}
\textit{Claim 1:} Let $a\in\mathbb{N}$ and $b\in\mathbb{N}$ satisfy $a+b=m$.
Then, there exist $g\in\mathbb{N}$, $x\in\mathbb{Z}$ and $y\in\mathbb{Z}$ such
that $g=ax+by$ and $g\mid a$ and $g\mid b$.
\end{statement}

Before we prove Claim 1, let us show a slightly weaker version of it, in which
we rename $a$ and $b$ as $u$ and $v$ and add the assumption that $u\geq v$:

\begin{statement}
\textit{Claim 2:} Let $u\in\mathbb{N}$ and $v\in\mathbb{N}$ satisfy $u+v=m$
and $u\geq v$. Then, there exist $g\in\mathbb{N}$, $x\in\mathbb{Z}$ and
$y\in\mathbb{Z}$ such that $g=ux+vy$ and $g\mid u$ and $g\mid v$.
\end{statement}

[\textit{Proof of Claim 2:} We are in one of the following two cases:

\textit{Case 1:} We have $v=0$.

\textit{Case 2:} We have $v\neq0$.

Let us first consider Case 1. In this case, we have $v=0$. Hence, $v=0=0u$, so
that $u\mid v$. Also, $u\cdot1+v\cdot0=u$. Thus, $u=u\cdot1+v\cdot0$ and
$u\mid u$ and $u\mid v$. Hence, there exist $g\in\mathbb{N}$, $x\in\mathbb{Z}$
and $y\in\mathbb{Z}$ such that $g=ux+vy$ and $g\mid u$ and $g\mid v$ (namely,
$g=u$, $x=1$ and $y=0$). Thus, Claim 2 is proven in Case 1.

Let us now consider Case 2. In this case, we have $v\neq0$. Hence, $v>0$
(since $v\in\mathbb{N}$). Thus, $u+v>u+0=u$, so that $u<u+v=m$. Hence,
(\ref{pf.thm.ind.bezout.IH}) (applied to $n=u$) yields that $\mathcal{A}%
\left(  u\right)  $ holds. In other words,%
\begin{equation}
\left(
\begin{array}
[c]{c}%
\text{if }a\in\mathbb{N}\text{ and }b\in\mathbb{N}\text{ satisfy }a+b=u\text{,
then there exist }g\in\mathbb{N}\text{, }x\in\mathbb{Z}\\
\text{and }y\in\mathbb{Z}\text{ such that }g=ax+by\text{ and }g\mid a\text{
and }g\mid b
\end{array}
\right)  \label{pf.thm.ind.bezout.C2.pf.2}%
\end{equation}
(because this is what the statement $\mathcal{A}\left(  u\right)  $ says).

Also, $u-v\in\mathbb{N}$ (since $u\geq v$) and $\left(  u-v\right)  +v=u$.
Hence, (\ref{pf.thm.ind.bezout.C2.pf.2}) (applied to $a=u-v$ and $b=v$) shows
that there exist $g\in\mathbb{N}$, $x\in\mathbb{Z}$ and $y\in\mathbb{Z}$ such
that $g=\left(  u-v\right)  x+vy$ and $g\mid u-v$ and $g\mid v$. Consider
these $g$, $x$ and $y$, and denote them by $g^{\prime}$, $x^{\prime}$ and
$y^{\prime}$. Thus, $g^{\prime}$ is an element of $\mathbb{N}$, and
$x^{\prime}$ and $y^{\prime}$ are elements of $\mathbb{Z}$ satisfying
$g^{\prime}=\left(  u-v\right)  x^{\prime}+vy^{\prime}$ and $g^{\prime}\mid
u-v$ and $g^{\prime}\mid v$.

Now, we have $g^{\prime}\mid u-v$; in other words, $u\equiv
v\operatorname{mod}g^{\prime}$. Also, $g^{\prime}\mid v$; in other words,
$v\equiv0\operatorname{mod}g^{\prime}$. Hence, $u\equiv v\equiv
0\operatorname{mod}g^{\prime}$, so that $u\equiv0\operatorname{mod}g^{\prime}%
$. In other words, $g^{\prime}\mid u$. Furthermore,
\[
g^{\prime}=\left(  u-v\right)  x^{\prime}+vy^{\prime}=ux^{\prime}-vx^{\prime
}+vy^{\prime}=ux^{\prime}+v\left(  y^{\prime}-x^{\prime}\right)  .
\]
Hence, there exist $g\in\mathbb{N}$, $x\in\mathbb{Z}$ and $y\in\mathbb{Z}$
such that $g=ux+vy$ and $g\mid u$ and $g\mid v$ (namely, $g=g^{\prime}$,
$x=x^{\prime}$ and $y=y^{\prime}-x^{\prime}$). Thus, Claim 2 is proven in Case 2.

We have now proven Claim 2 in each of the two Cases 1 and 2. Thus, Claim 2
always holds (since Cases 1 and 2 cover all possibilities).]

Now, we can prove Claim 1 as well:

[\textit{Proof of Claim 1:} We are in one of the following two cases:

\textit{Case 1:} We have $a\geq b$.

\textit{Case 2:} We have $a<b$.

Let us first consider Case 1. In this case, we have $a\geq b$. Hence, Claim 2
(applied to $u=a$ and $v=b$) shows that there exist $g\in\mathbb{N}$,
$x\in\mathbb{Z}$ and $y\in\mathbb{Z}$ such that $g=ax+by$ and $g\mid a$ and
$g\mid b$. Thus, Claim 1 is proven in Case 1.

Let us next consider Case 2. In this case, we have $a<b$. Hence, $a\leq b$, so
that $b\geq a$. Also, $b+a=a+b=m$. Hence, Claim 2 (applied to $u=b$ and $v=a$)
shows that there exist $g\in\mathbb{N}$, $x\in\mathbb{Z}$ and $y\in\mathbb{Z}$
such that $g=bx+ay$ and $g\mid b$ and $g\mid a$. Consider these $g$, $x$ and
$y$, and denote them by $g^{\prime}$, $x^{\prime}$ and $y^{\prime}$. Thus,
$g^{\prime}$ is an element of $\mathbb{N}$, and $x^{\prime}$ and $y^{\prime}$
are elements of $\mathbb{Z}$ satisfying $g^{\prime}=bx^{\prime}+ay^{\prime}$
and $g^{\prime}\mid b$ and $g^{\prime}\mid a$. Now, $g^{\prime}=bx^{\prime
}+ay^{\prime}=ay^{\prime}+bx^{\prime}$. Hence, there exist $g\in\mathbb{N}$,
$x\in\mathbb{Z}$ and $y\in\mathbb{Z}$ such that $g=ax+by$ and $g\mid a$ and
$g\mid b$ (namely, $g=g^{\prime}$, $x=y^{\prime}$ and $y=x^{\prime}$). Thus,
Claim 1 is proven in Case 2.

We have now proven Claim 1 in each of the two Cases 1 and 2. Thus, Claim 1
always holds (since Cases 1 and 2 cover all possibilities).]

But $\mathcal{A}\left(  m\right)  $ is defined as the statement%
\[
\left(
\begin{array}
[c]{c}%
\text{if }a\in\mathbb{N}\text{ and }b\in\mathbb{N}\text{ satisfy }a+b=m\text{,
then there exist }g\in\mathbb{N}\text{, }x\in\mathbb{Z}\\
\text{and }y\in\mathbb{Z}\text{ such that }g=ax+by\text{ and }g\mid a\text{
and }g\mid b
\end{array}
\right)  .
\]
Thus, $\mathcal{A}\left(  m\right)  $ is precisely Claim 1. Hence,
$\mathcal{A}\left(  m\right)  $ holds (since Claim 1 holds). This completes
the induction step. Thus, we have proven by strong induction that
$\mathcal{A}\left(  n\right)  $ holds for all $n\in\mathbb{N}$. In other
words, the statement (\ref{pf.thm.ind.bezout.An}) holds for all $n\in
\mathbb{N}$ (since this statement is precisely $\mathcal{A}\left(  n\right)  $).

Now, let $a\in\mathbb{N}$ and $b\in\mathbb{N}$. Then, $a+b\in\mathbb{N}$.
Hence, we can apply (\ref{pf.thm.ind.bezout.An}) to $n=a+b$ (since $a+b=a+b$).
We thus conclude that there exist $g\in\mathbb{N}$, $x\in\mathbb{Z}$ and
$y\in\mathbb{Z}$ such that $g=ax+by$ and $g\mid a$ and $g\mid b$. This proves
Theorem \ref{thm.ind.bezout}.
\end{proof}

\subsubsection{Conventions for writing proofs by strong induction on derived
quantities}

Let us take a closer look at the proof we just gave. The statement
$\mathcal{A}\left(  n\right)  $ that we defined was unsurprising: It simply
says that Theorem \ref{thm.ind.bezout} holds under the condition that $a+b=n$.
Thus, by introducing $\mathcal{A}\left(  n\right)  $, we have
\textquotedblleft sliced\textquotedblright\ Theorem \ref{thm.ind.bezout} into
a sequence of statements $\mathcal{A}\left(  0\right)  ,\mathcal{A}\left(
1\right)  ,\mathcal{A}\left(  2\right)  ,\ldots$, which then allowed us to
prove these statements by strong induction on $n$ even though no
\textquotedblleft$n$\textquotedblright\ appeared in Theorem
\ref{thm.ind.bezout} itself. This strong induction can be simply called a
\textquotedblleft strong induction on $a+b$\textquotedblright. More generally:

\begin{convention}
\label{conv.ind.SIPder}Let $\mathcal{B}$ be a logical statement that involves
some variables $v_{1},v_{2},v_{3},\ldots$. (For example, $\mathcal{B}$ can be
the statement of Theorem \ref{thm.ind.bezout}; then, these variables are $a$
and $b$.)

Let $g\in\mathbb{Z}$. (This $g$ has nothing to do with the $g$ from Theorem
\ref{thm.ind.bezout}.)

Let $q$ be some expression (involving the variables $v_{1},v_{2},v_{3},\ldots$
or some of them) that has the property that whenever the variables
$v_{1},v_{2},v_{3},\ldots$ satisfy the assumptions of $\mathcal{B}$, the
expression $q$ evaluates to some element of $\mathbb{Z}_{\geq g}$. (For
example, if $\mathcal{B}$ is the statement of Theorem \ref{thm.ind.bezout} and
$g=0$, then $q$ can be the expression $a+b$, because $a+b\in\mathbb{N}%
=\mathbb{Z}_{\geq0}$ whenever $a$ and $b$ are as in Theorem
\ref{thm.ind.bezout}.)

Assume that you want to prove the statement $\mathcal{B}$. Then, you can
proceed as follows: For each $n\in\mathbb{Z}_{\geq g}$, define $\mathcal{A}%
\left(  n\right)  $ to be the statement saying that\footnotemark%
\[
\left(  \text{the statement }\mathcal{B}\text{ holds under the condition that
}q=n\right)  .
\]
Then, prove $\mathcal{A}\left(  n\right)  $ by strong induction on $n$
starting at $g$. Thus:

\begin{itemize}
\item The \textit{induction step} consists in fixing $m\in\mathbb{Z}_{\geq g}%
$, and showing that if
\begin{equation}
\left(  \mathcal{A}\left(  n\right)  \text{ holds for every }n\in
\mathbb{Z}_{\geq g}\text{ satisfying }n<m\right)  ,
\label{eq.conv.ind.SIPder.IH1}%
\end{equation}
then
\begin{equation}
\left(  \mathcal{A}\left(  m\right)  \text{ holds}\right)  .
\label{eq.conv.ind.SIPder.IG1}%
\end{equation}
In other words, it consists in fixing $m\in\mathbb{Z}_{\geq g}$, and showing
that if
\begin{equation}
\left(  \text{the statement }\mathcal{B}\text{ holds under the condition that
}q<m\right)  , \label{eq.conv.ind.SIPder.IH2}%
\end{equation}
then
\begin{equation}
\left(  \text{the statement }\mathcal{B}\text{ holds under the condition that
}q=m\right)  . \label{eq.conv.ind.SIPder.IG2}%
\end{equation}
(Indeed, the previous two sentences are equivalent, because of the logical
equivalences%
\begin{align*}
&  \ \left(  \mathcal{A}\left(  n\right)  \text{ holds for every }%
n\in\mathbb{Z}_{\geq g}\text{ satisfying }n<m\right) \\
\Longleftrightarrow\  &  \left(
\begin{array}
[c]{c}%
\left(  \text{the statement }\mathcal{B}\text{ holds under the condition that
}q=n\right) \\
\text{holds for every }n\in\mathbb{Z}_{\geq g}\text{ satisfying }n<m
\end{array}
\right) \\
&  \ \ \ \ \ \ \ \ \ \ \left(
\begin{array}
[c]{c}%
\text{since the statement }\mathcal{A}\left(  n\right)  \text{ is defined
as}\\
\left(  \text{the statement }\mathcal{B}\text{ holds under the condition that
}q=n\right)
\end{array}
\right) \\
\Longleftrightarrow\  &  \left(  \text{the statement }\mathcal{B}\text{ holds
under the condition that }q<m\right)
\end{align*}
and%
\begin{align*}
&  \left(  \mathcal{A}\left(  m\right)  \text{ holds}\right) \\
\Longleftrightarrow\  &  \left(  \text{the statement }\mathcal{B}\text{ holds
under the condition that }q=m\right) \\
&  \ \ \ \ \ \ \ \ \ \ \left(
\begin{array}
[c]{c}%
\text{since the statement }\mathcal{A}\left(  m\right)  \text{ is defined
as}\\
\left(  \text{the statement }\mathcal{B}\text{ holds under the condition that
}q=m\right)
\end{array}
\right)  .
\end{align*}
)

In practice, this induction step will usually be organized as follows: We fix
$m\in\mathbb{Z}_{\geq g}$, then we assume that the statement $\mathcal{B}$
holds under the condition that $q<m$ (this is the induction hypothesis), and
then we prove that the statement $\mathcal{B}$ holds under the condition that
$q=m$.
\end{itemize}

Once this induction proof is finished, it immediately follows that the
statement $\mathcal{B}$ always holds (because the induction proof has shown
that, whatever $n\in\mathbb{Z}_{\geq g}$ is, the statement $\mathcal{B}$ holds
under the condition that $q=n$).

This strategy of proof is called \textquotedblleft strong induction on
$q$\textquotedblright\ (or \textquotedblleft strong induction over
$q$\textquotedblright). Once you have specified what $q$ is, you don't need to
explicitly define $\mathcal{A}\left(  n\right)  $, nor do you ever need to
mention $n$.
\end{convention}

\footnotetext{We assume that no variable named \textquotedblleft%
$n$\textquotedblright\ appears in the statement $\mathcal{B}$; otherwise, we
need a different letter for our new variable in order to avoid confusion.}%
Using this convention, we can rewrite our above proof of Theorem
\ref{thm.ind.bezout} as follows (remembering once again that $\mathbb{Z}%
_{\geq0}=\mathbb{N}$):

\begin{proof}
[Proof of Theorem \ref{thm.ind.bezout} (second version).]Let us prove Theorem
\ref{thm.ind.bezout} by strong induction on $a+b$ starting at $0$:

\textit{Induction step:} Let $m\in\mathbb{N}$. Assume that Theorem
\ref{thm.ind.bezout} holds under the condition that $a+b<m$. We must then show
that Theorem \ref{thm.ind.bezout} holds under the condition that $a+b=m$. This
is tantamount to proving the following claim:

\begin{statement}
\textit{Claim 1:} Let $a\in\mathbb{N}$ and $b\in\mathbb{N}$ satisfy $a+b=m$.
Then, there exist $g\in\mathbb{N}$, $x\in\mathbb{Z}$ and $y\in\mathbb{Z}$ such
that $g=ax+by$ and $g\mid a$ and $g\mid b$.
\end{statement}

Before we prove Claim 1, let us show a slightly weaker version of it, in which
we rename $a$ and $b$ as $u$ and $v$ and add the assumption that $u\geq v$:

\begin{statement}
\textit{Claim 2:} Let $u\in\mathbb{N}$ and $v\in\mathbb{N}$ satisfy $u+v=m$
and $u\geq v$. Then, there exist $g\in\mathbb{N}$, $x\in\mathbb{Z}$ and
$y\in\mathbb{Z}$ such that $g=ux+vy$ and $g\mid u$ and $g\mid v$.
\end{statement}

[\textit{Proof of Claim 2:} We are in one of the following two cases:

\textit{Case 1:} We have $v=0$.

\textit{Case 2:} We have $v\neq0$.

Let us first consider Case 1. In this case, we have $v=0$. Hence, $v=0=0u$, so
that $u\mid v$. Also, $u\cdot1+v\cdot0=u$. Thus, $u=u\cdot1+v\cdot0$ and
$u\mid u$ and $u\mid v$. Hence, there exist $g\in\mathbb{N}$, $x\in\mathbb{Z}$
and $y\in\mathbb{Z}$ such that $g=ux+vy$ and $g\mid u$ and $g\mid v$ (namely,
$g=u$, $x=1$ and $y=0$). Thus, Claim 2 is proven in Case 1.

Let us now consider Case 2. In this case, we have $v\neq0$. Hence, $v>0$
(since $v\in\mathbb{N}$). Thus, $u+v>u+0=u$, so that $u<u+v=m$. Also,
$u-v\in\mathbb{N}$ (since $u\geq v$) and $\left(  u-v\right)  +v=u$.

But we assumed that Theorem \ref{thm.ind.bezout} holds under the condition
that $a+b<m$. Thus, we can apply Theorem \ref{thm.ind.bezout} to $a=u-v$ and
$b=v$ (since $u-v\in\mathbb{N}$ and $\left(  u-v\right)  +v=u<m$). We thus
conclude that there exist $g\in\mathbb{N}$, $x\in\mathbb{Z}$ and
$y\in\mathbb{Z}$ such that $g=\left(  u-v\right)  x+vy$ and $g\mid u-v$ and
$g\mid v$. Consider these $g$, $x$ and $y$, and denote them by $g^{\prime}$,
$x^{\prime}$ and $y^{\prime}$. Thus, $g^{\prime}$ is an element of
$\mathbb{N}$, and $x^{\prime}$ and $y^{\prime}$ are elements of $\mathbb{Z}$
satisfying $g^{\prime}=\left(  u-v\right)  x^{\prime}+vy^{\prime}$ and
$g^{\prime}\mid u-v$ and $g^{\prime}\mid v$.

Now, we have $g^{\prime}\mid u-v$; in other words, $u\equiv
v\operatorname{mod}g^{\prime}$. Also, $g^{\prime}\mid v$; in other words,
$v\equiv0\operatorname{mod}g^{\prime}$. Hence, $u\equiv v\equiv
0\operatorname{mod}g^{\prime}$, so that $u\equiv0\operatorname{mod}g^{\prime}%
$. In other words, $g^{\prime}\mid u$. Furthermore,%
\[
g^{\prime}=\left(  u-v\right)  x^{\prime}+vy^{\prime}=ux^{\prime}-vx^{\prime
}+vy^{\prime}=ux^{\prime}+v\left(  y^{\prime}-x^{\prime}\right)  .
\]
Hence, there exist $g\in\mathbb{N}$, $x\in\mathbb{Z}$ and $y\in\mathbb{Z}$
such that $g=ux+vy$ and $g\mid u$ and $g\mid v$ (namely, $g=g^{\prime}$,
$x=x^{\prime}$ and $y=y^{\prime}-x^{\prime}$). Thus, Claim 2 is proven in Case 2.

We have now proven Claim 2 in each of the two Cases 1 and 2. Thus, Claim 2
always holds (since Cases 1 and 2 cover all possibilities).]

Now, we can prove Claim 1 as well:

[\textit{Proof of Claim 1:} Claim 1 can be derived from Claim 2 in the same
way as we derived it in the first version of the proof above. We shall not
repeat this argument, since it just applies verbatim.]

But Claim 1 is simply saying that Theorem \ref{thm.ind.bezout} holds under the
condition that $a+b=m$. Thus, by proving Claim 1, we have shown that Theorem
\ref{thm.ind.bezout} holds under the condition that $a+b=m$. This completes
the induction step. Thus, Theorem \ref{thm.ind.bezout} is proven by strong induction.
\end{proof}

\subsection{\label{sect.ind.interval}Induction in an interval}

\subsubsection{The induction principle for intervals}

The induction principles we have seen so far were tailored towards proving
statements whose variables range over infinite sets such as $\mathbb{N}$ and
$\mathbb{Z}_{\geq g}$. Sometimes, one instead wants to do an induction on a
variable that ranges over a finite interval, such as $\left\{  g,g+1,\ldots
,h\right\}  $ for some integers $g$ and $h$. We shall next state an induction
principle tailored to such situations. First, we make an important convention:

\begin{convention}
If $g$ and $h$ are two integers such that $g>h$, then the set $\left\{
g,g+1,\ldots,h\right\}  $ is understood to be the empty set.
\end{convention}

Thus, for example, $\left\{  2,3,\ldots,1\right\}  =\varnothing$ and $\left\{
2,3,\ldots,0\right\}  =\varnothing$ and $\left\{  5,6,\ldots,-100\right\}
=\varnothing$. (But $\left\{  5,6,\ldots,5\right\}  =\left\{  5\right\}  $ and
$\left\{  5,6,\ldots,6\right\}  =\left\{  5,6\right\}  $.)

We now state our induction principle for intervals:

\begin{theorem}
\label{thm.ind.IPgh}Let $g\in\mathbb{Z}$ and $h\in\mathbb{Z}$. For each
$n\in\left\{  g,g+1,\ldots,h\right\}  $, let $\mathcal{A}\left(  n\right)  $
be a logical statement.

Assume the following:

\begin{statement}
\textit{Assumption 1:} If $g\leq h$, then the statement $\mathcal{A}\left(
g\right)  $ holds.
\end{statement}

\begin{statement}
\textit{Assumption 2:} If $m\in\left\{  g,g+1,\ldots,h-1\right\}  $ is such
that $\mathcal{A}\left(  m\right)  $ holds, then $\mathcal{A}\left(
m+1\right)  $ also holds.
\end{statement}

Then, $\mathcal{A}\left(  n\right)  $ holds for each $n\in\left\{
g,g+1,\ldots,h\right\}  $.
\end{theorem}

Theorem \ref{thm.ind.IPgh} is, in a sense, the closest one can get to Theorem
\ref{thm.ind.IPg} when having only finitely many statements $\mathcal{A}%
\left(  g\right)  ,\mathcal{A}\left(  g+1\right)  ,\ldots,\mathcal{A}\left(
h\right)  $ instead of an infinite sequence of statements $\mathcal{A}\left(
g\right)  ,\mathcal{A}\left(  g+1\right)  ,\mathcal{A}\left(  g+2\right)
,\ldots$. It is easy to derive Theorem \ref{thm.ind.IPgh} from Corollary
\ref{cor.ind.IPg.renamed}:

\begin{proof}
[Proof of Theorem \ref{thm.ind.IPgh}.]For each $n\in\mathbb{Z}_{\geq g}$, we
define $\mathcal{B}\left(  n\right)  $ to be the logical statement%
\[
\left(  \text{if }n\in\left\{  g,g+1,\ldots,h\right\}  \text{, then
}\mathcal{A}\left(  n\right)  \text{ holds}\right)  .
\]

Now, let us consider the Assumptions A and B from Corollary
\ref{cor.ind.IPg.renamed}. We claim that both of these assumptions are satisfied.

Assumption 1 says that if $g\leq h$, then the statement $\mathcal{A}\left(
g\right)  $ holds. Thus, $\mathcal{B}\left(  g\right)  $
holds\footnote{\textit{Proof.} Assume that $g\in\left\{  g,g+1,\ldots
,h\right\}  $. Thus, $g\leq h$. But Assumption 1 says that if $g\leq h$, then
the statement $\mathcal{A}\left(  g\right)  $ holds. Hence, the statement
$\mathcal{A}\left(  g\right)  $ holds (since $g\leq h$).
\par
Now, forget that we assumed that $g\in\left\{  g,g+1,\ldots,h\right\}  $. We
thus have proven that if $g\in\left\{  g,g+1,\ldots,h\right\}  $, then
$\mathcal{A}\left(  g\right)  $ holds. In other words, $\mathcal{B}\left(
g\right)  $ holds (because the statement $\mathcal{B}\left(  g\right)  $ is
defined as $\left(  \text{if }g\in\left\{  g,g+1,\ldots,h\right\}  \text{,
then }\mathcal{A}\left(  g\right)  \text{ holds}\right)  $). Qed.}. In other
words, Assumption A is satisfied.

Next, we shall prove that Assumption B is satisfied. Indeed, let
$p\in\mathbb{Z}_{\geq g}$ be such that $\mathcal{B}\left(  p\right)  $ holds.
We shall now show that $\mathcal{B}\left(  p+1\right)  $ also holds.

Indeed, assume that $p+1\in\left\{  g,g+1,\ldots,h\right\}  $. Thus, $p+1\leq
h$, so that $p\leq p+1\leq h$. Combining this with $p\geq g$ (since
$p\in\mathbb{Z}_{\geq g}$), we conclude that $p\in\left\{  g,g+1,\ldots
,h\right\}  $ (since $p$ is an integer). But we have assumed that
$\mathcal{B}\left(  p\right)  $ holds. In other words,%
\[
\text{if }p\in\left\{  g,g+1,\ldots,h\right\}  \text{, then }\mathcal{A}%
\left(  p\right)  \text{ holds}%
\]
(because the statement $\mathcal{B}\left(  p\right)  $ is defined as $\left(
\text{if }p\in\left\{  g,g+1,\ldots,h\right\}  \text{, then }\mathcal{A}%
\left(  p\right)  \text{ holds}\right)  $). Thus, $\mathcal{A}\left(
p\right)  $ holds (since we have $p\in\left\{  g,g+1,\ldots,h\right\}  $).
Also, from $p+1\leq h$, we obtain $p\leq h-1$. Combining this with $p\geq g$,
we find $p\in\left\{  g,g+1,\ldots,h-1\right\}  $. Thus, we know that
$p\in\left\{  g,g+1,\ldots,h-1\right\}  $ is such that $\mathcal{A}\left(
p\right)  $ holds. Hence, Assumption 2 (applied to $m=p$) shows that
$\mathcal{A}\left(  p+1\right)  $ also holds.

Now, forget that we assumed that $p+1\in\left\{  g,g+1,\ldots,h\right\}  $. We
thus have proven that if $p+1\in\left\{  g,g+1,\ldots,h\right\}  $, then
$\mathcal{A}\left(  p+1\right)  $ holds. In other words, $\mathcal{B}\left(
p+1\right)  $ holds (since the statement $\mathcal{B}\left(  p+1\right)  $ is
defined as \newline$\left(  \text{if }p+1\in\left\{  g,g+1,\ldots,h\right\}
\text{, then }\mathcal{A}\left(  p+1\right)  \text{ holds}\right)  $).

Now, forget that we fixed $p$. We thus have proven that if $p\in
\mathbb{Z}_{\geq g}$ is such that $\mathcal{B}\left(  p\right)  $ holds, then
$\mathcal{B}\left(  p+1\right)  $ also holds. In other words, Assumption B is satisfied.

We now know that both Assumption A and Assumption B are satisfied. Hence,
Corollary \ref{cor.ind.IPg.renamed} shows that
\begin{equation}
\mathcal{B}\left(  n\right)  \text{ holds for each }n\in\mathbb{Z}_{\geq g}.
\label{pf.thm.ind.IPgh.at}%
\end{equation}

Now, let $n\in\left\{  g,g+1,\ldots,h\right\}  $. Thus, $n\geq g$, so that
$n\in\mathbb{Z}_{\geq g}$. Hence, (\ref{pf.thm.ind.IPgh.at}) shows that
$\mathcal{B}\left(  n\right)  $ holds. In other words,
\[
\text{if }n\in\left\{  g,g+1,\ldots,h\right\}  \text{, then }\mathcal{A}%
\left(  n\right)  \text{ holds}%
\]
(since the statement $\mathcal{B}\left(  n\right)  $ was defined as $\left(
\text{if }n\in\left\{  g,g+1,\ldots,h\right\}  \text{, then }\mathcal{A}%
\left(  n\right)  \text{ holds}\right)  $). Thus, $\mathcal{A}\left(
n\right)  $ holds (since we have $n\in\left\{  g,g+1,\ldots,h\right\}  $).

Now, forget that we fixed $n$. We thus have shown that $\mathcal{A}\left(
n\right)  $ holds for each $n\in\left\{  g,g+1,\ldots,h\right\}  $. This
proves Theorem \ref{thm.ind.IPgh}.
\end{proof}

Theorem \ref{thm.ind.IPgh} is called the \textit{principle of induction
starting at }$g$\textit{ and ending at }$h$, and proofs that use it are
usually called \textit{proofs by induction} or \textit{induction proofs}. As
with all the other induction principles seen so far, we don't usually
explicitly cite Theorem \ref{thm.ind.IPgh}, but instead say certain words that
signal that it is being applied and that (ideally) also indicate what integers
$g$ and $h$ and what statements $\mathcal{A}\left(  n\right)  $ it is being
applied to\footnote{We will explain this in Convention \ref{conv.ind.IPghlang}
below.}. However, we shall reference it explicitly in our very first example
of the use of Theorem \ref{thm.ind.IPgh}:

\begin{proposition}
\label{prop.ind.dc}Let $g$ and $h$ be integers such that $g\leq h$. Let
$b_{g},b_{g+1},\ldots,b_{h}$ be any $h-g+1$ nonzero integers. Assume that
$b_{g}\geq0$. Assume further that%
\begin{equation}
\left\vert b_{i+1}-b_{i}\right\vert \leq1\ \ \ \ \ \ \ \ \ \ \text{for every
}i\in\left\{  g,g+1,\ldots,h-1\right\}  . \label{eq.prop.ind.dc.ass}%
\end{equation}
Then, $b_{n}>0$ for each $n\in\left\{  g,g+1,\ldots,h\right\}  $.
\end{proposition}

Proposition \ref{prop.ind.dc} is often called the \textquotedblleft%
\textit{discrete intermediate value theorem}\textquotedblright\ or the
\textquotedblleft\textit{discrete continuity principle}\textquotedblright. Its
intuitive meaning is that if a finite list of nonzero integers starts with a
nonnegative integer, and every further entry of this list differs from its
preceding entry by at most $1$, then all entries of this list must be
positive. An example of such a list is $\left(
2,3,3,2,3,4,4,3,2,3,2,3,2,1\right)  $. Notice that Proposition
\ref{prop.ind.dc} is, again, rather obvious from an intuitive perspective: It
just says that it isn't possible to go from a nonnegative integer to a
negative integer by steps of $1$ without ever stepping at $0$. The rigorous
proof of Proposition \ref{prop.ind.dc} is not much harder -- but because it is
a statement about elements of $\left\{  g,g+1,\ldots,h\right\}  $, it
naturally relies on Theorem \ref{thm.ind.IPgh}:

\begin{proof}
[Proof of Proposition \ref{prop.ind.dc}.]For each $n\in\left\{  g,g+1,\ldots
,h\right\}  $, we let $\mathcal{A}\left(  n\right)  $ be the statement
$\left(  b_{n}>0\right)  $.

Our next goal is to prove the statement $\mathcal{A}\left(  n\right)  $ for
each $n\in\left\{  g,g+1,\ldots,h\right\}  $.

All the $h-g+1$ integers $b_{g},b_{g+1},\ldots,b_{h}$ are nonzero (by
assumption). Thus, in particular, $b_{g}$ is nonzero. In other words,
$b_{g}\neq0$. Combining this with $b_{g}\geq0$, we obtain $b_{g}>0$. In other
words, the statement $\mathcal{A}\left(  g\right)  $ holds (since this
statement $\mathcal{A}\left(  g\right)  $ is defined to be $\left(
b_{g}>0\right)  $). Hence,
\begin{equation}
\text{if }g\leq h\text{, then the statement }\mathcal{A}\left(  g\right)
\text{ holds.} \label{pf.prop.ind.dc.base}%
\end{equation}

Now, we claim that
\begin{equation}
\text{if }m\in\left\{  g,g+1,\ldots,h-1\right\}  \text{ is such that
}\mathcal{A}\left(  m\right)  \text{ holds, then }\mathcal{A}\left(
m+1\right)  \text{ also holds.} \label{pf.prop.ind.dc.step}%
\end{equation}

[\textit{Proof of (\ref{pf.prop.ind.dc.step}):} Let $m\in\left\{
g,g+1,\ldots,h-1\right\}  $ be such that $\mathcal{A}\left(  m\right)  $
holds. We must show that $\mathcal{A}\left(  m+1\right)  $ also holds.

We have assumed that $\mathcal{A}\left(  m\right)  $ holds. In other words,
$b_{m}>0$ holds (since $\mathcal{A}\left(  m\right)  $ is defined to be the
statement $\left(  b_{m}>0\right)  $). Now, (\ref{eq.prop.ind.dc.ass})
(applied to $i=m$) yields $\left\vert b_{m+1}-b_{m}\right\vert \leq1$. But it
is well-known (and easy to see) that every integer $x$ satisfies
$-x\leq\left\vert x\right\vert $. Applying this to $x=b_{m+1}-b_{m}$, we
obtain $-\left(  b_{m+1}-b_{m}\right)  \leq\left\vert b_{m+1}-b_{m}\right\vert
\leq1$. In other words, $1\geq-\left(  b_{m+1}-b_{m}\right)  =b_{m}-b_{m+1}$.
In other words, $1+b_{m+1}\geq b_{m}$. Hence, $1+b_{m+1}\geq b_{m}>0$, so that
$1+b_{m+1}\geq1$ (since $1+b_{m+1}$ is an integer). In other words,
$b_{m+1}\geq0$.

But all the $h-g+1$ integers $b_{g},b_{g+1},\ldots,b_{h}$ are nonzero (by
assumption). Thus, in particular, $b_{m+1}$ is nonzero. In other words,
$b_{m+1}\neq0$. Combining this with $b_{m+1}\geq0$, we obtain $b_{m+1}>0$. But
this is precisely the statement $\mathcal{A}\left(  m+1\right)  $ (because
$\mathcal{A}\left(  m+1\right)  $ is defined to be the statement $\left(
b_{m+1}>0\right)  $). Thus, the statement $\mathcal{A}\left(  m+1\right)  $ holds.

Now, forget that we fixed $m$. We thus have shown that if $m\in\left\{
g,g+1,\ldots,h-1\right\}  $ is such that $\mathcal{A}\left(  m\right)  $
holds, then $\mathcal{A}\left(  m+1\right)  $ also holds. This proves
(\ref{pf.prop.ind.dc.step}).]

Now, both assumptions of Theorem \ref{thm.ind.IPgh} are satisfied (indeed,
Assumption 1 holds because of (\ref{pf.prop.ind.dc.base}), whereas Assumption
2 holds because of (\ref{pf.prop.ind.dc.step})). Thus, Theorem
\ref{thm.ind.IPgh} shows that $\mathcal{A}\left(  n\right)  $ holds for each
$n\in\left\{  g,g+1,\ldots,h\right\}  $. In other words, $b_{n}>0$ holds for
each $n\in\left\{  g,g+1,\ldots,h\right\}  $ (since $\mathcal{A}\left(
n\right)  $ is the statement $\left(  b_{n}>0\right)  $). This proves
Proposition \ref{prop.ind.dc}.
\end{proof}

\subsubsection{Conventions for writing induction proofs in intervals}

Next, we shall introduce some standard language that is commonly used in
proofs by induction starting at $g$ and ending at $h$. This language closely
imitates the one we use for proofs by standard induction:

\begin{convention}
\label{conv.ind.IPghlang}Let $g\in\mathbb{Z}$ and $h\in\mathbb{Z}$. For each
$n\in\left\{  g,g+1,\ldots,h\right\}  $, let $\mathcal{A}\left(  n\right)  $
be a logical statement. Assume that you want to prove that $\mathcal{A}\left(
n\right)  $ holds for each $n\in\left\{  g,g+1,\ldots,h\right\}  $.

Theorem \ref{thm.ind.IPgh} offers the following strategy for proving this:
First show that Assumption 1 of Theorem \ref{thm.ind.IPgh} is satisfied; then,
show that Assumption 2 of Theorem \ref{thm.ind.IPgh} is satisfied; then,
Theorem \ref{thm.ind.IPgh} automatically completes your proof.

A proof that follows this strategy is called a \textit{proof by induction on
}$n$ (or \textit{proof by induction over }$n$) \textit{starting at }$g$
\textit{and ending at }$h$ or (less precisely) an \textit{inductive proof}.
Most of the time, the words \textquotedblleft starting at $g$ and ending at
$h$\textquotedblright\ are omitted, since they merely repeat what is clear
from the context anyway: For example, if you make a claim about all integers
$n\in\left\{  3,4,5,6\right\}  $, and you say that you are proving it by
induction on $n$, it is clear that you are using induction on $n$ starting at
$3$ and ending at $6$.

The proof that Assumption 1 is satisfied is called the \textit{induction base}
(or \textit{base case}) of the proof. The proof that Assumption 2 is satisfied
is called the \textit{induction step} of the proof.

In order to prove that Assumption 2 is satisfied, you will usually want to fix
an $m\in\left\{  g,g+1,\ldots,h-1\right\}  $ such that $\mathcal{A}\left(
m\right)  $ holds, and then prove that $\mathcal{A}\left(  m+1\right)  $
holds. In other words, you will usually want to fix $m\in\left\{
g,g+1,\ldots,h-1\right\}  $, assume that $\mathcal{A}\left(  m\right)  $
holds, and then prove that $\mathcal{A}\left(  m+1\right)  $ holds. When doing
so, it is common to refer to the assumption that $\mathcal{A}\left(  m\right)
$ holds as the \textit{induction hypothesis} (or \textit{induction assumption}).
\end{convention}

Unsurprisingly, this language parallels the language introduced in Convention
\ref{conv.ind.IP0lang} and in Convention \ref{conv.ind.IPglang}.

Again, we can shorten our inductive proofs by omitting some sentences that
convey no information. In particular, we can leave out the explicit definition
of the statement $\mathcal{A}\left(  n\right)  $ when this statement is
precisely the claim that we are proving (without the \textquotedblleft for
each $n\in\left\{  g,g+1,\ldots,h\right\}  $\textquotedblright\ part).
Furthermore, it is common to leave the \textquotedblleft If $g\leq
h$\textquotedblright\ part of Assumption 1 unsaid (i.e., to pretend that
Assumption 1 simply says that $\mathcal{A}\left(  g\right)  $ holds). Strictly
speaking, this is somewhat imprecise, since $\mathcal{A}\left(  g\right)  $ is
not defined when $g>h$; but of course, the whole claim that is being proven is
moot anyway when $g>h$ (because there exist no $n\in\left\{  g,g+1,\ldots
,h\right\}  $ in this case), so this imprecision doesn't matter.

Thus, we can rewrite our above proof of Proposition \ref{prop.ind.dc} as follows:

\begin{proof}
[Proof of Proposition \ref{prop.ind.dc} (second version).]We claim that%
\begin{equation}
b_{n}>0 \label{pf.prop.ind.dc.2nd.claim}%
\end{equation}
for each $n\in\left\{  g,g+1,\ldots,h\right\}  $.

Indeed, we shall prove (\ref{pf.prop.ind.dc.2nd.claim}) by induction on $n$:

\textit{Induction base:} All the $h-g+1$ integers $b_{g},b_{g+1},\ldots,b_{h}$
are nonzero (by assumption). Thus, in particular, $b_{g}$ is nonzero. In other
words, $b_{g}\neq0$. Combining this with $b_{g}\geq0$, we obtain $b_{g}>0$. In
other words, (\ref{pf.prop.ind.dc.2nd.claim}) holds for $n=g$. This completes
the induction base.

\textit{Induction step:} Let $m\in\left\{  g,g+1,\ldots,h-1\right\}  $. Assume
that (\ref{pf.prop.ind.dc.2nd.claim}) holds for $n=m$. We must show that
(\ref{pf.prop.ind.dc.2nd.claim}) also holds for $n=m+1$.

We have assumed that (\ref{pf.prop.ind.dc.2nd.claim}) holds for $n=m$. In
other words, $b_{m}>0$. Now, (\ref{eq.prop.ind.dc.ass}) (applied to $i=m$)
yields $\left\vert b_{m+1}-b_{m}\right\vert \leq1$. But it is well-known (and
easy to see) that every integer $x$ satisfies $-x\leq\left\vert x\right\vert
$. Applying this to $x=b_{m+1}-b_{m}$, we obtain $-\left(  b_{m+1}%
-b_{m}\right)  \leq\left\vert b_{m+1}-b_{m}\right\vert \leq1$. In other words,
$1\geq-\left(  b_{m+1}-b_{m}\right)  =b_{m}-b_{m+1}$. In other words,
$1+b_{m+1}\geq b_{m}$. Hence, $1+b_{m+1}\geq b_{m}>0$, so that $1+b_{m+1}%
\geq1$ (since $1+b_{m+1}$ is an integer). In other words, $b_{m+1}\geq0$.

But all the $h-g+1$ integers $b_{g},b_{g+1},\ldots,b_{h}$ are nonzero (by
assumption). Thus, in particular, $b_{m+1}$ is nonzero. In other words,
$b_{m+1}\neq0$. Combining this with $b_{m+1}\geq0$, we obtain $b_{m+1}>0$. In
other words, (\ref{pf.prop.ind.dc.2nd.claim}) holds for $n=m+1$. This
completes the induction step. Thus, (\ref{pf.prop.ind.dc.2nd.claim}) is proven
by induction. This proves Proposition \ref{prop.ind.dc}.
\end{proof}

\subsection{\label{sect.ind.strong-interval}Strong induction in an interval}

\subsubsection{The strong induction principle for intervals}

We shall next state yet another induction principle -- one that combines the
idea of strong induction (as in Theorem \ref{thm.ind.SIP}) with the idea of
working inside an interval $\left\{  g,g+1,\ldots,h\right\}  $ (as in Theorem
\ref{thm.ind.IPgh}):

\begin{theorem}
\label{thm.ind.SIPgh}Let $g\in\mathbb{Z}$ and $h\in\mathbb{Z}$. For each
$n\in\left\{  g,g+1,\ldots,h\right\}  $, let $\mathcal{A}\left(  n\right)  $
be a logical statement.

Assume the following:

\begin{statement}
\textit{Assumption 1:} If $m\in\left\{  g,g+1,\ldots,h\right\}  $ is such that%
\[
\left(  \mathcal{A}\left(  n\right)  \text{ holds for every }n\in\left\{
g,g+1,\ldots,h\right\}  \text{ satisfying }n<m\right)  ,
\]
then $\mathcal{A}\left(  m\right)  $ holds.
\end{statement}

Then, $\mathcal{A}\left(  n\right)  $ holds for each $n\in\left\{
g,g+1,\ldots,h\right\}  $.
\end{theorem}

Our proof of Theorem \ref{thm.ind.SIPgh} will be similar to the proof of
Theorem \ref{thm.ind.IPgh}, except that we shall be using Theorem
\ref{thm.ind.SIP} instead of Corollary \ref{cor.ind.IPg.renamed}. Or, to be
more precise, we shall be using the following restatement of Theorem
\ref{thm.ind.SIP}:

\begin{corollary}
\label{cor.ind.SIP.renamed}Let $g\in\mathbb{Z}$. For each $n\in\mathbb{Z}%
_{\geq g}$, let $\mathcal{B}\left(  n\right)  $ be a logical statement.

Assume the following:

\begin{statement}
\textit{Assumption A:} If $p\in\mathbb{Z}_{\geq g}$ is such that%
\[
\left(  \mathcal{B}\left(  n\right)  \text{ holds for every }n\in
\mathbb{Z}_{\geq g}\text{ satisfying }n<p\right)  ,
\]
then $\mathcal{B}\left(  p\right)  $ holds.
\end{statement}

Then, $\mathcal{B}\left(  n\right)  $ holds for each $n\in\mathbb{Z}_{\geq g}$.
\end{corollary}

\begin{proof}
[Proof of Corollary \ref{cor.ind.SIP.renamed}.]Corollary
\ref{cor.ind.SIP.renamed} is exactly Theorem \ref{thm.ind.SIP}, except that
some names have been changed:

\begin{itemize}
\item The statements $\mathcal{A}\left(  n\right)  $ have been renamed as
$\mathcal{B}\left(  n\right)  $.

\item Assumption 1 has been renamed as Assumption A.

\item The variable $m$ in Assumption A has been renamed as $p$.
\end{itemize}

Thus, Corollary \ref{cor.ind.SIP.renamed} holds (since Theorem
\ref{thm.ind.SIP} holds).
\end{proof}

We can now prove Theorem \ref{thm.ind.SIPgh}:

\begin{proof}
[Proof of Theorem \ref{thm.ind.SIPgh}.]For each $n\in\mathbb{Z}_{\geq g}$, we
define $\mathcal{B}\left(  n\right)  $ to be the logical statement%
\[
\left(  \text{if }n\in\left\{  g,g+1,\ldots,h\right\}  \text{, then
}\mathcal{A}\left(  n\right)  \text{ holds}\right)  .
\]

Now, let us consider the Assumption A from Corollary \ref{cor.ind.SIP.renamed}%
. We claim that this assumption is satisfied.

Indeed, let $p\in\mathbb{Z}_{\geq g}$ be such that
\begin{equation}
\left(  \mathcal{B}\left(  n\right)  \text{ holds for every }n\in
\mathbb{Z}_{\geq g}\text{ satisfying }n<p\right)  .
\label{pf.thm.ind.SIPgh.AA}%
\end{equation}
We shall now show that $\mathcal{B}\left(  p\right)  $ holds.

Indeed, assume that $p\in\left\{  g,g+1,\ldots,h\right\}  $. Thus, $p\leq h$.

Now, let $n\in\left\{  g,g+1,\ldots,h\right\}  $ be such that $n<p$. Then,
$n\in\left\{  g,g+1,\ldots,h\right\}  \subseteq\left\{  g,g+1,g+2,\ldots
\right\}  =\mathbb{Z}_{\geq g}$ and $n<p$. Hence, (\ref{pf.thm.ind.SIPgh.AA})
shows that $\mathcal{B}\left(  n\right)  $ holds. In other words, $\left(
\text{if }n\in\left\{  g,g+1,\ldots,h\right\}  \text{, then }\mathcal{A}%
\left(  n\right)  \text{ holds}\right)  $ (because the statement
$\mathcal{B}\left(  n\right)  $ is defined as $\left(  \text{if }n\in\left\{
g,g+1,\ldots,h\right\}  \text{, then }\mathcal{A}\left(  n\right)  \text{
holds}\right)  $). Therefore, $\mathcal{A}\left(  n\right)  $ holds (since we
know that $n\in\left\{  g,g+1,\ldots,h\right\}  $).

Now, forget that we fixed $n$. We thus have proven that
\[
\left(  \mathcal{A}\left(  n\right)  \text{ holds for every }n\in\left\{
g,g+1,\ldots,h\right\}  \text{ satisfying }n<p\right)  .
\]
Hence, Assumption 1 (applied to $m=p$) yields that $\mathcal{A}\left(
p\right)  $ holds.

Now, forget that we assumed that $p\in\left\{  g,g+1,\ldots,h\right\}  $. We
thus have proven that
\[
\left(  \text{if }p\in\left\{  g,g+1,\ldots,h\right\}  \text{, then
}\mathcal{A}\left(  p\right)  \text{ holds}\right)  .
\]
In other words, $\mathcal{B}\left(  p\right)  $ holds (since the statement
$\mathcal{B}\left(  p\right)  $ was defined as \newline$\left(  \text{if }%
p\in\left\{  g,g+1,\ldots,h\right\}  \text{, then }\mathcal{A}\left(
p\right)  \text{ holds}\right)  $).

Now, forget that we fixed $p$. We thus have shown that if $p\in\mathbb{Z}%
_{\geq g}$ is such that%
\[
\left(  \mathcal{B}\left(  n\right)  \text{ holds for every }n\in
\mathbb{Z}_{\geq g}\text{ satisfying }n<p\right)  ,
\]
then $\mathcal{B}\left(  p\right)  $ holds. In other words, Assumption A is satisfied.

Hence, Corollary \ref{cor.ind.SIP.renamed} shows that%
\begin{equation}
\mathcal{B}\left(  n\right)  \text{ holds for each }n\in\mathbb{Z}_{\geq g}.
\label{pf.thm.ind.SIPgh.at}%
\end{equation}

Now, let $n\in\left\{  g,g+1,\ldots,h\right\}  $. Thus, $n\geq g$, so that
$n\in\mathbb{Z}_{\geq g}$. Hence, (\ref{pf.thm.ind.SIPgh.at}) shows that
$\mathcal{B}\left(  n\right)  $ holds. In other words,
\[
\text{if }n\in\left\{  g,g+1,\ldots,h\right\}  \text{, then }\mathcal{A}%
\left(  n\right)  \text{ holds}%
\]
(since the statement $\mathcal{B}\left(  n\right)  $ was defined as $\left(
\text{if }n\in\left\{  g,g+1,\ldots,h\right\}  \text{, then }\mathcal{A}%
\left(  n\right)  \text{ holds}\right)  $). Thus, $\mathcal{A}\left(
n\right)  $ holds (since we have $n\in\left\{  g,g+1,\ldots,h\right\}  $).

Now, forget that we fixed $n$. We thus have shown that $\mathcal{A}\left(
n\right)  $ holds for each $n\in\left\{  g,g+1,\ldots,h\right\}  $. This
proves Theorem \ref{thm.ind.SIPgh}.
\end{proof}

Theorem \ref{thm.ind.SIPgh} is called the \textit{principle of strong
induction starting at }$g$\textit{ and ending at }$h$, and proofs that use it
are usually called \textit{proofs by strong induction}. Once again, we usually
don't explicitly cite Theorem \ref{thm.ind.SIPgh} in such proofs, and we
usually don't say explicitly what $g$ and $h$ are and what the statements
$\mathcal{A}\left(  n\right)  $ are when it is clear from the context. But (as
with all the other induction principles considered so far) we shall be
explicit about all these details in our first example:

\begin{proposition}
\label{prop.ind.dcs}Let $g$ and $h$ be integers such that $g\leq h$. Let
$b_{g},b_{g+1},\ldots,b_{h}$ be any $h-g+1$ nonzero integers. Assume that
$b_{g}\geq0$. Assume that for each $p\in\left\{  g+1,g+2,\ldots,h\right\}  $,
\begin{equation}
\text{there exists some }j\in\left\{  g,g+1,\ldots,p-1\right\}  \text{ such
that }b_{p}\geq b_{j}-1. \label{eq.prop.ind.dcs.ass}%
\end{equation}
(Of course, the $j$ can depend on $p$.) Then, $b_{n}>0$ for each $n\in\left\{
g,g+1,\ldots,h\right\}  $.
\end{proposition}

Proposition \ref{prop.ind.dcs} is a more general (although less intuitive)
version of Proposition \ref{prop.ind.dc}; indeed, it is easy to see that the
condition (\ref{eq.prop.ind.dc.ass}) is stronger than the condition
(\ref{eq.prop.ind.dcs.ass}) (when required for all $p\in\left\{
g+1,g+2,\ldots,h\right\}  $).

\begin{example}
For this example, set $g=3$ and $h=7$. Then, if we set $\left(  b_{3}%
,b_{4},b_{5},b_{6},b_{7}\right)  =\left(  4,5,3,4,2\right)  $, then the
condition (\ref{eq.prop.ind.dcs.ass}) holds for all $p\in\left\{
g+1,g+2,\ldots,h\right\}  $. (For example, it holds for $p=5$, since
$b_{5}=3\geq4-1=b_{1}-1$ and $1\in\left\{  g,g+1,\ldots,5-1\right\}  $.) On
the other hand, if we set $\left(  b_{3},b_{4},b_{5},b_{6},b_{7}\right)
=\left(  4,5,2,4,3\right)  $, then this condition does not hold (indeed, it
fails for $p=5$, since $b_{5}=2$ is neither $\geq4-1$ nor $\geq5-1$).
\end{example}

Let us now prove Proposition \ref{prop.ind.dcs} using Theorem
\ref{thm.ind.SIPgh}:

\begin{proof}
[Proof of Proposition \ref{prop.ind.dcs}.]For each $n\in\left\{
g,g+1,\ldots,h\right\}  $, we let $\mathcal{A}\left(  n\right)  $ be the
statement $\left(  b_{n}>0\right)  $.

Our next goal is to prove the statement $\mathcal{A}\left(  n\right)  $ for
each $n\in\left\{  g,g+1,\ldots,h\right\}  $.

All the $h-g+1$ integers $b_{g},b_{g+1},\ldots,b_{h}$ are nonzero (by
assumption). Thus, in particular, $b_{g}$ is nonzero. In other words,
$b_{g}\neq0$. Combining this with $b_{g}\geq0$, we obtain $b_{g}>0$. In other
words, the statement $\mathcal{A}\left(  g\right)  $ holds (since this
statement $\mathcal{A}\left(  g\right)  $ is defined to be $\left(
b_{g}>0\right)  $).

Now, we make the following claim:

\begin{statement}
\textit{Claim 1:} If $m\in\left\{  g,g+1,\ldots,h\right\}  $ is such that%
\[
\left(  \mathcal{A}\left(  n\right)  \text{ holds for every }n\in\left\{
g,g+1,\ldots,h\right\}  \text{ satisfying }n<m\right)  ,
\]
then $\mathcal{A}\left(  m\right)  $ holds.
\end{statement}

[\textit{Proof of Claim 1:} Let $m\in\left\{  g,g+1,\ldots,h\right\}  $ be
such that
\begin{equation}
\left(  \mathcal{A}\left(  n\right)  \text{ holds for every }n\in\left\{
g,g+1,\ldots,h\right\}  \text{ satisfying }n<m\right)  .
\label{pf.prop.ind.dcs.3}%
\end{equation}
We must show that $\mathcal{A}\left(  m\right)  $ holds.

If $m=g$, then this follows from the fact that $\mathcal{A}\left(  g\right)  $
holds. Thus, for the rest of the proof of Claim 1, we WLOG assume that we
don't have $m=g$. Hence, $m\neq g$. Combining this with $m\in\left\{
g,g+1,\ldots,h\right\}  $, we obtain $m\in\left\{  g,g+1,\ldots,h\right\}
\setminus\left\{  g\right\}  \subseteq\left\{  g+1,g+2,\ldots,h\right\}  $.
Hence, (\ref{eq.prop.ind.dcs.ass}) (applied to $p=m$) shows that there exists
some $j\in\left\{  g,g+1,\ldots,m-1\right\}  $ such that $b_{m}\geq b_{j}-1$.
Consider this $j$. From $m\in\left\{  g+1,g+2,\ldots,h\right\}  $, we obtain
$m\leq h$.

From $j\in\left\{  g,g+1,\ldots,m-1\right\}  $, we obtain $j\leq m-1<m$. Also,
\newline$j\in\left\{  g,g+1,\ldots,m-1\right\}  \subseteq\left\{
g,g+1,\ldots,h\right\}  $ (since $m-1\leq m\leq h$). Thus,
(\ref{pf.prop.ind.dcs.3}) (applied to $n=j$) yields that $\mathcal{A}\left(
j\right)  $ holds. In other words, $b_{j}>0$ holds (since $\mathcal{A}\left(
j\right)  $ is defined to be the statement $\left(  b_{j}>0\right)  $). Thus,
$b_{j}\geq1$ (since $b_{j}$ is an integer), so that $b_{j}-1\geq0$. But recall
that $b_{m}\geq b_{j}-1\geq0$.

But all the $h-g+1$ integers $b_{g},b_{g+1},\ldots,b_{h}$ are nonzero (by
assumption). Thus, in particular, $b_{m}$ is nonzero. In other words,
$b_{m}\neq0$. Combining this with $b_{m}\geq0$, we obtain $b_{m}>0$. But this
is precisely the statement $\mathcal{A}\left(  m\right)  $ (because
$\mathcal{A}\left(  m\right)  $ is defined to be the statement $\left(
b_{m}>0\right)  $). Thus, the statement $\mathcal{A}\left(  m\right)  $ holds.
This completes the proof of Claim 1.]

Claim 1 says that Assumption 1 of Theorem \ref{thm.ind.SIPgh} is satisfied.
Thus, Theorem \ref{thm.ind.SIPgh} shows that $\mathcal{A}\left(  n\right)  $
holds for each $n\in\left\{  g,g+1,\ldots,h\right\}  $. In other words,
$b_{n}>0$ holds for each $n\in\left\{  g,g+1,\ldots,h\right\}  $ (since
$\mathcal{A}\left(  n\right)  $ is the statement $\left(  b_{n}>0\right)  $).
This proves Proposition \ref{prop.ind.dcs}.
\end{proof}

\subsubsection{Conventions for writing strong induction proofs in intervals}

Next, we shall introduce some standard language that is commonly used in
proofs by strong induction starting at $g$ and ending at $h$. This language
closely imitates the one we use for proofs by \textquotedblleft
usual\textquotedblright\ strong induction:

\begin{convention}
\label{conv.ind.SIPghlang}Let $g\in\mathbb{Z}$ and $h\in\mathbb{Z}$. For each
$n\in\left\{  g,g+1,\ldots,h\right\}  $, let $\mathcal{A}\left(  n\right)  $
be a logical statement. Assume that you want to prove that $\mathcal{A}\left(
n\right)  $ holds for each $n\in\left\{  g,g+1,\ldots,h\right\}  $.

Theorem \ref{thm.ind.SIPgh} offers the following strategy for proving this:
Show that Assumption 1 of Theorem \ref{thm.ind.SIPgh} is satisfied; then,
Theorem \ref{thm.ind.SIPgh} automatically completes your proof.

A proof that follows this strategy is called a \textit{proof by strong
induction on }$n$ \textit{starting at }$g$ \textit{and ending at }$h$. Most of
the time, the words \textquotedblleft starting at $g$ and ending at
$h$\textquotedblright\ are omitted. The proof that Assumption 1 is satisfied
is called the \textit{induction step} of the proof. This kind of proof has no
\textquotedblleft induction base\textquotedblright.

In order to prove that Assumption 1 is satisfied, you will usually want to fix
an $m\in\left\{  g,g+1,\ldots,h\right\}  $ such that
\begin{equation}
\left(  \mathcal{A}\left(  n\right)  \text{ holds for every }n\in\left\{
g,g+1,\ldots,h\right\}  \text{ satisfying }n<m\right)  ,
\label{eq.conv.ind.SIPghlang.IH}%
\end{equation}
and then prove that $\mathcal{A}\left(  m\right)  $ holds. In other words, you
will usually want to fix $m\in\left\{  g,g+1,\ldots,h\right\}  $, assume that
(\ref{eq.conv.ind.SIPghlang.IH}) holds, and then prove that $\mathcal{A}%
\left(  m\right)  $ holds. When doing so, it is common to refer to the
assumption that (\ref{eq.conv.ind.SIPghlang.IH}) holds as the
\textit{induction hypothesis} (or \textit{induction assumption}).
\end{convention}

Unsurprisingly, this language parallels the language introduced in Convention
\ref{conv.ind.SIPlang}.

As before, proofs using strong induction can be shortened by leaving out some
uninformative prose. In particular, the explicit definition of the statement
$\mathcal{A}\left(  n\right)  $ can often be omitted when this statement is
precisely the claim that we are proving (without the \textquotedblleft for
each $n\in\left\{  g,g+1,\ldots,h\right\}  $\textquotedblright\ part). The
values of $g$ and $h$ can also be inferred from the statement of the claim, so
they don't need to be specified explicitly. And once again, we don't need to
write \textquotedblleft\textit{Induction step:}\textquotedblright, since our
strong induction has no induction base.

This leads to the following abridged version of our above proof of Proposition
\ref{prop.ind.dcs}:

\begin{proof}
[Proof of Proposition \ref{prop.ind.dcs} (second version).]We claim that%
\begin{equation}
b_{n}>0 \label{pf.prop.ind.dcs.2nd.claim}%
\end{equation}
for each $n\in\left\{  g,g+1,\ldots,h\right\}  $.

Indeed, we shall prove (\ref{pf.prop.ind.dcs.2nd.claim}) by strong induction
on $n$:

Let $m\in\left\{  g,g+1,\ldots,h\right\}  $. Assume that
(\ref{pf.prop.ind.dcs.2nd.claim}) holds for every $n\in\left\{  g,g+1,\ldots
,h\right\}  $ satisfying $n<m$. We must show that
(\ref{pf.prop.ind.dcs.2nd.claim}) also holds for $n=m$. In other words, we
must show that $b_{m}>0$.

All the $h-g+1$ integers $b_{g},b_{g+1},\ldots,b_{h}$ are nonzero (by
assumption). Thus, in particular, $b_{g}$ is nonzero. In other words,
$b_{g}\neq0$. Combining this with $b_{g}\geq0$, we obtain $b_{g}>0$.

We have assumed that (\ref{pf.prop.ind.dcs.2nd.claim}) holds for every
$n\in\left\{  g,g+1,\ldots,h\right\}  $ satisfying $n<m$. In other words, we
have%
\begin{equation}
b_{n}>0\text{ for every }n\in\left\{  g,g+1,\ldots,h\right\}  \text{
satisfying }n<m. \label{pf.prop.ind.dcs.2nd.claim.IH}%
\end{equation}

Recall that we must prove that $b_{m}>0$. If $m=g$, then this follows from
$b_{g}>0$. Thus, for the rest of this induction step, we WLOG assume that we
don't have $m=g$. Hence, $m\neq g$. Combining this with $m\in\left\{
g,g+1,\ldots,h\right\}  $, we obtain $m\in\left\{  g,g+1,\ldots,h\right\}
\setminus\left\{  g\right\}  \subseteq\left\{  g+1,g+2,\ldots,h\right\}  $.
Hence, (\ref{eq.prop.ind.dcs.ass}) (applied to $p=m$) shows that there exists
some $j\in\left\{  g,g+1,\ldots,m-1\right\}  $ such that $b_{m}\geq b_{j}-1$.
Consider this $j$. From $m\in\left\{  g+1,g+2,\ldots,h\right\}  $, we obtain
$m\leq h$.

From $j\in\left\{  g,g+1,\ldots,m-1\right\}  $, we obtain $j\leq m-1<m$. Also,
\newline$j\in\left\{  g,g+1,\ldots,m-1\right\}  \subseteq\left\{
g,g+1,\ldots,h\right\}  $ (since $m-1\leq m\leq h$). Thus,
(\ref{pf.prop.ind.dcs.2nd.claim.IH}) (applied to $n=j$) yields that $b_{j}>0$.
Thus, $b_{j}\geq1$ (since $b_{j}$ is an integer), so that $b_{j}-1\geq0$. But
recall that $b_{m}\geq b_{j}-1\geq0$.

But all the $h-g+1$ integers $b_{g},b_{g+1},\ldots,b_{h}$ are nonzero (by
assumption). Thus, in particular, $b_{m}$ is nonzero. In other words,
$b_{m}\neq0$. Combining this with $b_{m}\geq0$, we obtain $b_{m}>0$.

Thus, we have proven that $b_{m}>0$. In other words,
(\ref{pf.prop.ind.dcs.2nd.claim}) holds for $n=m$. This completes the
induction step. Thus, (\ref{pf.prop.ind.dcs.2nd.claim}) is proven by strong
induction. This proves Proposition \ref{prop.ind.dcs}.
\end{proof}

\subsection{\label{sect.ind.gen-ass}General associativity for composition of
maps}

\subsubsection{Associativity of map composition}

Recall that if $f:X\rightarrow Y$ and $g:Y\rightarrow Z$ are two maps, then
the \textit{composition} $g\circ f$ of the maps $g$ and $f$ is defined to be
the map%
\[
X\rightarrow Z,\ x\mapsto g\left(  f\left(  x\right)  \right)  .
\]

Now, if we have four sets $X$, $Y$, $Z$ and $W$ and three maps $c:X\rightarrow
Y$, $b:Y\rightarrow Z$ and $a:Z\rightarrow W$, then we can build two possible
compositions that use all three of these maps: namely, the two compositions
$\left(  a\circ b\right)  \circ c$ and $a\circ\left(  b\circ c\right)  $. It
turns out that these two compositions are the same map:\footnote{Of course,
when some of the four sets $X$, $Y$, $Z$ and $W$ are equal, then more
compositions can be built: For example, if $Y=Z=W$, then we can also build the
composition $\left(  b\circ a\right)  \circ c$ or the composition $\left(
\left(  b\circ b\right)  \circ a\right)  \circ c$. But these compositions are
not the same map as the two that we previously constructed.}

\begin{proposition}
\label{prop.ind.gen-ass-maps.fgh}Let $X$, $Y$, $Z$ and $W$ be four sets. Let
$c:X\rightarrow Y$, $b:Y\rightarrow Z$ and $a:Z\rightarrow W$ be three maps.
Then,%
\[
\left(  a\circ b\right)  \circ c=a\circ\left(  b\circ c\right)  .
\]

\end{proposition}

Proposition \ref{prop.ind.gen-ass-maps.fgh} is called the
\textit{associativity of map composition}, and is proven straightforwardly:

\begin{proof}
[Proof of Proposition \ref{prop.ind.gen-ass-maps.fgh}.]Let $x\in X$. Then, the
definition of $b\circ c$ yields $\left(  b\circ c\right)  \left(  x\right)
=b\left(  c\left(  x\right)  \right)  $. But the definition of $\left(  a\circ
b\right)  \circ c$ yields%
\[
\left(  \left(  a\circ b\right)  \circ c\right)  \left(  x\right)  =\left(
a\circ b\right)  \left(  c\left(  x\right)  \right)  =a\left(  b\left(
c\left(  x\right)  \right)  \right)  \ \ \ \ \ \ \ \ \ \ \left(  \text{by the
definition of }a\circ b\right)  .
\]
On the other hand, the definition of $a\circ\left(  b\circ c\right)  $ yields%
\[
\left(  a\circ\left(  b\circ c\right)  \right)  \left(  x\right)  =a\left(
\underbrace{\left(  b\circ c\right)  \left(  x\right)  }_{=b\left(  c\left(
x\right)  \right)  }\right)  =a\left(  b\left(  c\left(  x\right)  \right)
\right)  .
\]
Comparing these two equalities, we obtain $\left(  \left(  a\circ b\right)
\circ c\right)  \left(  x\right)  =\left(  a\circ\left(  b\circ c\right)
\right)  \left(  x\right)  $.

Now, forget that we fixed $x$. We thus have shown that
\[
\left(  \left(  a\circ b\right)  \circ c\right)  \left(  x\right)  =\left(
a\circ\left(  b\circ c\right)  \right)  \left(  x\right)
\ \ \ \ \ \ \ \ \ \ \text{for each }x\in X.
\]
In other words, $\left(  a\circ b\right)  \circ c=a\circ\left(  b\circ
c\right)  $. This proves Proposition \ref{prop.ind.gen-ass-maps.fgh}.
\end{proof}

\subsubsection{Composing more than $3$ maps: exploration}

Proposition \ref{prop.ind.gen-ass-maps.fgh} can be restated as follows: If
$a$, $b$ and $c$ are three maps such that the compositions $a\circ b$ and
$b\circ c$ are well-defined, then $\left(  a\circ b\right)  \circ
c=a\circ\left(  b\circ c\right)  $. This allows us to write \textquotedblleft%
$a\circ b\circ c$\textquotedblright\ for each of the compositions $\left(
a\circ b\right)  \circ c$ and $a\circ\left(  b\circ c\right)  $ without having
to disambiguate this expression by means of parentheses. It is natural to ask
whether we can do the same thing for more than three maps. For example, let us
consider four maps $a$, $b$, $c$ and $d$ for which the compositions $a\circ
b$, $b\circ c$ and $c\circ d$ are well-defined:

\begin{example}
\label{exa.ind.gen-ass-maps.abcd}Let $X$, $Y$, $Z$, $W$ and $U$ be five sets.
Let $d:X\rightarrow Y$, $c:Y\rightarrow Z$, $b:Z\rightarrow W$ and
$a:W\rightarrow U$ be four maps. Then, there we can construct five
compositions that use all four of these maps; these five compositions are%
\begin{align}
&  \left(  \left(  a\circ b\right)  \circ c\right)  \circ
d,\ \ \ \ \ \ \ \ \ \ \left(  a\circ\left(  b\circ c\right)  \right)  \circ
d,\ \ \ \ \ \ \ \ \ \ \left(  a\circ b\right)  \circ\left(  c\circ d\right)
,\label{eq.exa.ind.gen-ass-maps.abcd.cp1}\\
&  a\circ\left(  \left(  b\circ c\right)  \circ d\right)
,\ \ \ \ \ \ \ \ \ \ a\circ\left(  b\circ\left(  c\circ d\right)  \right)  .
\label{eq.exa.ind.gen-ass-maps.abcd.cp2}%
\end{align}
It turns out that these five compositions are all the same map. Indeed, this
follows by combining the following observations:

\begin{itemize}
\item We have $\left(  \left(  a\circ b\right)  \circ c\right)  \circ
d=\left(  a\circ\left(  b\circ c\right)  \right)  \circ d$ (since Proposition
\ref{prop.ind.gen-ass-maps.fgh} yields $\left(  a\circ b\right)  \circ
c=a\circ\left(  b\circ c\right)  $).

\item We have $a\circ\left(  \left(  b\circ c\right)  \circ d\right)
=a\circ\left(  b\circ\left(  c\circ d\right)  \right)  $ (since Proposition
\ref{prop.ind.gen-ass-maps.fgh} yields $\left(  b\circ c\right)  \circ
d=b\circ\left(  c\circ d\right)  $).

\item We have $\left(  a\circ\left(  b\circ c\right)  \right)  \circ
d=a\circ\left(  \left(  b\circ c\right)  \circ d\right)  $ (by Proposition
\ref{prop.ind.gen-ass-maps.fgh}, applied to $W$, $U$, $b\circ c$ and $d$
instead of $Z$, $W$, $b$ and $c$).

\item We have $\left(  \left(  a\circ b\right)  \circ c\right)  \circ
d=\left(  a\circ b\right)  \circ\left(  c\circ d\right)  $ (by Proposition
\ref{prop.ind.gen-ass-maps.fgh}, applied to $U$, $a\circ b$, $c$ and $d$
instead of $W$, $a$, $b$ and $c$).
\end{itemize}

Hence, all five compositions are equal. Thus, we can write \textquotedblleft%
$a\circ b\circ c\circ d$\textquotedblright\ for each of these five
compositions, again dropping the parentheses.

We shall refer to the five compositions listed in
(\ref{eq.exa.ind.gen-ass-maps.abcd.cp1}) and
(\ref{eq.exa.ind.gen-ass-maps.abcd.cp2}) as the \textquotedblleft complete
parenthesizations of $a\circ b\circ c\circ d$\textquotedblright. Here, the
word \textquotedblleft parenthesization\textquotedblright\ means a way to put
parentheses into the expression \textquotedblleft$a\circ b\circ c\circ
d$\textquotedblright, whereas the word \textquotedblleft
complete\textquotedblright\ means that these parentheses unambiguously
determine which two maps any given $\circ$ sign is composing. (For example,
the parenthesization \textquotedblleft$\left(  a\circ b\circ c\right)  \circ
d$\textquotedblright\ is not complete, because the first $\circ$ sign in it
could be either composing $a$ with $b$ or composing $a$ with $b\circ c$. But
the parenthesization \textquotedblleft$\left(  \left(  a\circ b\right)  \circ
c\right)  \circ d$\textquotedblright\ is complete, because its first $\circ$
sign composes $a$ and $b$, whereas its second $\circ$ sign composes $a\circ b$
with $c$, and finally its third $\circ$ sign composes $\left(  a\circ
b\right)  \circ c$ with $d$.)

Thus, we have seen that all five complete parenthesizations of $a\circ b\circ
c\circ d$ are the same map.
\end{example}

What happens if we compose more than four maps? Clearly, the more maps we
have, the more complete parenthesizations can be constructed. We have good
reasons to suspect that these parenthesizations will all be the same map (so
we can again drop the parentheses); but if we try to prove it in the ad-hoc
way we did in Example \ref{exa.ind.gen-ass-maps.abcd}, then we have more and
more work to do the more maps we are composing. Clearly, if we want to prove
our suspicion for arbitrarily many maps, we need a more general approach.

\subsubsection{Formalizing general associativity}

So let us make a general statement; but first, let us formally define the
notion of a \textquotedblleft complete parenthesization\textquotedblright:

\begin{definition}
\label{def.ind.gen-ass-maps.cp}Let $n$ be a positive integer. Let $X_{1}%
,X_{2},\ldots,X_{n+1}$ be $n+1$ sets. For each $i\in\left\{  1,2,\ldots
,n\right\}  $, let $f_{i}:X_{i}\rightarrow X_{i+1}$ be a map. Then, we want to
define the notion of a \textit{complete parenthesization} of $f_{n}\circ
f_{n-1}\circ\cdots\circ f_{1}$. We define this notion by recursion on $n$ as follows:

\begin{itemize}
\item For $n=1$, there is only one complete parenthesization of $f_{n}\circ
f_{n-1}\circ\cdots\circ f_{1}$, and this is simply the map $f_{1}%
:X_{1}\rightarrow X_{2}$.

\item If $n>1$, then the complete parenthesizations of $f_{n}\circ
f_{n-1}\circ\cdots\circ f_{1}$ are all the maps of the form $\alpha\circ\beta
$, where

\begin{itemize}
\item $k$ is some element of $\left\{  1,2,\ldots,n-1\right\}  $;

\item $\alpha$ is a complete parenthesization of $f_{n}\circ f_{n-1}%
\circ\cdots\circ f_{k+1}$;

\item $\beta$ is a complete parenthesization of $f_{k}\circ f_{k-1}\circ
\cdots\circ f_{1}$.
\end{itemize}
\end{itemize}
\end{definition}

\begin{example}
Let us see what this definition yields for small values of $n$:

\begin{itemize}
\item For $n=1$, the only complete parenthesization of $f_{1}$ is $f_{1}$.

\item For $n=2$, the only complete parenthesization of $f_{2}\circ f_{1}$ is
the composition $f_{2}\circ f_{1}$ (because here, the only possible values of
$k$, $\alpha$ and $\beta$ are $1$, $f_{2}$ and $f_{1}$, respectively).

\item For $n=3$, the complete parenthesizations of $f_{3}\circ f_{2}\circ
f_{1}$ are the two compositions $\left(  f_{3}\circ f_{2}\right)  \circ f_{1}$
and $f_{3}\circ\left(  f_{2}\circ f_{1}\right)  $ (because here, the only
possible values of $k$ are $1$ and $2$, and each value of $k$ uniquely
determines $\alpha$ and $\beta$). Proposition \ref{prop.ind.gen-ass-maps.fgh}
shows that they are equal (as maps).

\item For $n=4$, the complete parenthesizations of $f_{4}\circ f_{3}\circ
f_{2}\circ f_{1}$ are the five compositions%
\begin{align*}
&  \left(  \left(  f_{4}\circ f_{3}\right)  \circ f_{2}\right)  \circ
f_{1},\ \ \ \ \ \ \ \ \ \ \left(  f_{4}\circ\left(  f_{3}\circ f_{2}\right)
\right)  \circ f_{1},\ \ \ \ \ \ \ \ \ \ \left(  f_{4}\circ f_{3}\right)
\circ\left(  f_{2}\circ f_{1}\right)  ,\\
&  f_{4}\circ\left(  \left(  f_{3}\circ f_{2}\right)  \circ f_{1}\right)
,\ \ \ \ \ \ \ \ \ \ f_{4}\circ\left(  f_{3}\circ\left(  f_{2}\circ
f_{1}\right)  \right)  .
\end{align*}
(These are exactly the five compositions listed in
(\ref{eq.exa.ind.gen-ass-maps.abcd.cp1}) and
(\ref{eq.exa.ind.gen-ass-maps.abcd.cp2}), except that the maps $d,c,b,a$ are
now called $f_{1},f_{2},f_{3},f_{4}$.) We have seen in Example
\ref{exa.ind.gen-ass-maps.abcd} that these five compositions are equal as maps.

\item For $n=5$, the complete parenthesizations of $f_{5}\circ f_{4}\circ
f_{3}\circ f_{2}\circ f_{1}$ are $14$ compositions, one of which is $\left(
f_{5}\circ f_{4}\right)  \circ\left(  f_{3}\circ\left(  f_{2}\circ
f_{1}\right)  \right)  $. Again, it is laborious but not difficult to check
that all the $14$ compositions are equal as maps.
\end{itemize}
\end{example}

Now, we want to prove the following general statement:

\begin{theorem}
\label{thm.ind.gen-ass-maps.cp}Let $n$ be a positive integer. Let $X_{1}%
,X_{2},\ldots,X_{n+1}$ be $n+1$ sets. For each $i\in\left\{  1,2,\ldots
,n\right\}  $, let $f_{i}:X_{i}\rightarrow X_{i+1}$ be a map. Then, all
complete parenthesizations of $f_{n}\circ f_{n-1}\circ\cdots\circ f_{1}$ are
the same map (from $X_{1}$ to $X_{n+1}$).
\end{theorem}

Theorem \ref{thm.ind.gen-ass-maps.cp} is sometimes called the \textit{general
associativity} theorem, and is often proved in the context of monoids (see,
e.g., \cite[Proposition 2.1.4]{Artin}); while the context is somewhat
different from ours, the proofs usually given still apply in ours.

\subsubsection{Defining the \textquotedblleft canonical\textquotedblright%
\ composition $C\left(  f_{n},f_{n-1},\ldots,f_{1}\right)  $}

We shall prove Theorem \ref{thm.ind.gen-ass-maps.cp} in a slightly indirect
way: We first define a \textit{specific} complete parenthesization of
$f_{n}\circ f_{n-1}\circ\cdots\circ f_{1}$, which we shall call $C\left(
f_{n},f_{n-1},\ldots,f_{1}\right)  $; then we will show that it satisfies
certain equalities (Proposition \ref{prop.ind.gen-ass-maps.Ceq}), and then
prove that every complete parenthesization of $f_{n}\circ f_{n-1}\circ
\cdots\circ f_{1}$ equals this map $C\left(  f_{n},f_{n-1},\ldots
,f_{1}\right)  $ (Proposition \ref{prop.ind.gen-ass-maps.Ceq-cp}). Each step
of this strategy will rely on induction.

We begin with the definition of $C\left(  f_{n},f_{n-1},\ldots,f_{1}\right)  $:

\begin{definition}
\label{def.ind.gen-ass-maps.C}Let $n$ be a positive integer. Let $X_{1}%
,X_{2},\ldots,X_{n+1}$ be $n+1$ sets. For each $i\in\left\{  1,2,\ldots
,n\right\}  $, let $f_{i}:X_{i}\rightarrow X_{i+1}$ be a map. Then, we want to
define a map $C\left(  f_{n},f_{n-1},\ldots,f_{1}\right)  :X_{1}\rightarrow
X_{n+1}$. We define this map by recursion on $n$ as follows:

\begin{itemize}
\item If $n=1$, then we define $C\left(  f_{n},f_{n-1},\ldots,f_{1}\right)  $
to be the map $f_{1}:X_{1}\rightarrow X_{2}$. (Note that in this case,
$C\left(  f_{n},f_{n-1},\ldots,f_{1}\right)  =C\left(  f_{1}\right)  $,
because $\left(  f_{n},f_{n-1},\ldots,f_{1}\right)  =\left(  f_{1}%
,f_{1-1},\ldots,f_{1}\right)  =\left(  f_{1}\right)  $.)

\item If $n>1$, then we define $C\left(  f_{n},f_{n-1},\ldots,f_{1}\right)
:X_{1}\rightarrow X_{n+1}$ by%
\begin{equation}
C\left(  f_{n},f_{n-1},\ldots,f_{1}\right)  =f_{n}\circ C\left(
f_{n-1},f_{n-2},\ldots,f_{1}\right)  . \label{eq.def.ind.gen-ass-maps.C.rec}%
\end{equation}

\end{itemize}
\end{definition}

\begin{example}
\label{exa.ind.gen-ass-maps.Cex}Consider the situation of Definition
\ref{def.ind.gen-ass-maps.C}.

\textbf{(a)} If $n=1$, then
\begin{equation}
C\left(  f_{1}\right)  =f_{1} \label{eq.exa.ind.gen-ass-maps.Cex.1}%
\end{equation}
(by the $n=1$ case of the definition).

\textbf{(b)} If $n=2$, then%
\begin{align}
C\left(  f_{2},f_{1}\right)   &  =f_{2}\circ\underbrace{C\left(  f_{1}\right)
}_{\substack{=f_{1}\\\text{(by (\ref{eq.exa.ind.gen-ass-maps.Cex.1}))}%
}}\ \ \ \ \ \ \ \ \ \ \left(  \text{by (\ref{eq.def.ind.gen-ass-maps.C.rec}),
applied to }n=2\right) \nonumber\\
&  =f_{2}\circ f_{1}. \label{eq.exa.ind.gen-ass-maps.Cex.2}%
\end{align}

\textbf{(c)} If $n=3$, then
\begin{align}
C\left(  f_{3},f_{2},f_{1}\right)   &  =f_{3}\circ\underbrace{C\left(
f_{2},f_{1}\right)  }_{\substack{=f_{2}\circ f_{1}\\\text{(by
(\ref{eq.exa.ind.gen-ass-maps.Cex.2}))}}}\ \ \ \ \ \ \ \ \ \ \left(  \text{by
(\ref{eq.def.ind.gen-ass-maps.C.rec}), applied to }n=3\right) \nonumber\\
&  =f_{3}\circ\left(  f_{2}\circ f_{1}\right)  .
\label{eq.exa.ind.gen-ass-maps.Cex.3}%
\end{align}

\textbf{(d)} If $n=4$, then%
\begin{align}
C\left(  f_{4},f_{3},f_{2},f_{1}\right)   &  =f_{4}\circ\underbrace{C\left(
f_{3},f_{2},f_{1}\right)  }_{\substack{=f_{3}\circ\left(  f_{2}\circ
f_{1}\right)  \\\text{(by (\ref{eq.exa.ind.gen-ass-maps.Cex.3}))}%
}}\ \ \ \ \ \ \ \ \ \ \left(  \text{by (\ref{eq.def.ind.gen-ass-maps.C.rec}),
applied to }n=4\right) \nonumber\\
&  =f_{4}\circ\left(  f_{3}\circ\left(  f_{2}\circ f_{1}\right)  \right)  .
\label{eq.exa.ind.gen-ass-maps.Cex.4}%
\end{align}

\textbf{(e)} For an arbitrary $n\geq1$, we can informally write $C\left(
f_{n},f_{n-1},\ldots,f_{1}\right)  $ as%
\[
C\left(  f_{n},f_{n-1},\ldots,f_{1}\right)  =f_{n}\circ\left(  f_{n-1}%
\circ\left(  f_{n-2}\circ\left(  \cdots\circ\left(  f_{2}\circ f_{1}\right)
\cdots\right)  \right)  \right)  .
\]
The right hand side of this equality is a complete parenthesization of
$f_{n}\circ f_{n-1}\circ\cdots\circ f_{1}$, where all the parentheses are
\textquotedblleft concentrated as far right as possible\textquotedblright%
\ (i.e., there is an opening parenthesis after each \textquotedblleft$\circ
$\textquotedblright\ sign except for the last one; and there are $n-2$ closing
parentheses at the end of the expression). This is merely a visual restatement
of the recursive definition of $C\left(  f_{n},f_{n-1},\ldots,f_{1}\right)  $
we gave above.
\end{example}

\subsubsection{The crucial property of $C\left(  f_{n},f_{n-1},\ldots
,f_{1}\right)  $}

The following proposition will be key to our proof of Theorem
\ref{thm.ind.gen-ass-maps.cp}:

\begin{proposition}
\label{prop.ind.gen-ass-maps.Ceq}Let $n$ be a positive integer. Let
$X_{1},X_{2},\ldots,X_{n+1}$ be $n+1$ sets. For each $i\in\left\{
1,2,\ldots,n\right\}  $, let $f_{i}:X_{i}\rightarrow X_{i+1}$ be a map. Then,%
\[
C\left(  f_{n},f_{n-1},\ldots,f_{1}\right)  =C\left(  f_{n},f_{n-1}%
,\ldots,f_{k+1}\right)  \circ C\left(  f_{k},f_{k-1},\ldots,f_{1}\right)
\]
for each $k\in\left\{  1,2,\ldots,n-1\right\}  $.
\end{proposition}

\begin{proof}
[Proof of Proposition \ref{prop.ind.gen-ass-maps.Ceq}.]Forget that we fixed
$n$, $X_{1},X_{2},\ldots,X_{n+1}$ and the maps $f_{i}$. We shall prove
Proposition \ref{prop.ind.gen-ass-maps.Ceq} by induction on $n$%
:\ \ \ \ \footnote{The induction principle that we are applying here is
Theorem \ref{thm.ind.IPg} with $g=1$ (since $\mathbb{Z}_{\geq1}$ is the set of
all positive integers).}

\textit{Induction base:} If $n=1$, then $\left\{  1,2,\ldots,n-1\right\}
=\left\{  1,2,\ldots,1-1\right\}  =\varnothing$. Hence, if $n=1$, then there
exists no $k\in\left\{  1,2,\ldots,n-1\right\}  $. Thus, if $n=1$, then
Proposition \ref{prop.ind.gen-ass-maps.Ceq} is vacuously true (since
Proposition \ref{prop.ind.gen-ass-maps.Ceq} has a \textquotedblleft for each
$k\in\left\{  1,2,\ldots,n-1\right\}  $\textquotedblright\ clause). This
completes the induction base.

\textit{Induction step:} Let $m\in\mathbb{Z}_{\geq1}$. Assume that Proposition
\ref{prop.ind.gen-ass-maps.Ceq} holds under the condition that $n=m$. We must
now prove that Proposition \ref{prop.ind.gen-ass-maps.Ceq} holds under the
condition that $n=m+1$. In other words, we must prove the following claim:

\begin{statement}
\textit{Claim 1:} Let $X_{1},X_{2},\ldots,X_{\left(  m+1\right)  +1}$ be
$\left(  m+1\right)  +1$ sets. For each $i\in\left\{  1,2,\ldots,m+1\right\}
$, let $f_{i}:X_{i}\rightarrow X_{i+1}$ be a map. Then,
\[
C\left(  f_{m+1},f_{\left(  m+1\right)  -1},\ldots,f_{1}\right)  =C\left(
f_{m+1},f_{\left(  m+1\right)  -1},\ldots,f_{k+1}\right)  \circ C\left(
f_{k},f_{k-1},\ldots,f_{1}\right)
\]
for each $k\in\left\{  1,2,\ldots,\left(  m+1\right)  -1\right\}  $.
\end{statement}

[\textit{Proof of Claim 1:} Let $k\in\left\{  1,2,\ldots,\left(  m+1\right)
-1\right\}  $. Thus, $k\in\left\{  1,2,\ldots,\left(  m+1\right)  -1\right\}
=\left\{  1,2,\ldots,m\right\}  $ (since $\left(  m+1\right)  -1=m$).

We know that $X_{1},X_{2},\ldots,X_{\left(  m+1\right)  +1}$ are $\left(
m+1\right)  +1$ sets. In other words, \newline$X_{1},X_{2},\ldots,X_{m+2}$ are
$m+2$ sets (since $\left(  m+1\right)  +1=m+2$). We have $m\in\mathbb{Z}%
_{\geq1}$, thus $m\geq1>0$; hence, $m+1>1$. Thus,
(\ref{eq.def.ind.gen-ass-maps.C.rec}) (applied to $n=m+1$) yields%
\begin{align}
C\left(  f_{m+1},f_{\left(  m+1\right)  -1},\ldots,f_{1}\right)   &
=f_{m+1}\circ C\left(  f_{\left(  m+1\right)  -1},f_{\left(  m+1\right)
-2},\ldots,f_{1}\right) \nonumber\\
&  =f_{m+1}\circ C\left(  f_{m},f_{m-1},\ldots,f_{1}\right)
\label{pf.prop.ind.gen-ass-maps.Ceq.c1.pf.1}%
\end{align}
(since $\left(  m+1\right)  -1=m$ and $\left(  m+1\right)  -2=m-1$).

But we are in one of the following two cases:

\textit{Case 1:} We have $k=m$.

\textit{Case 2:} We have $k\neq m$.

Let us first consider Case 1. In this case, we have $k=m$. Hence,%
\[
C\left(  f_{m+1},f_{\left(  m+1\right)  -1},\ldots,f_{k+1}\right)  =C\left(
f_{m+1},f_{\left(  m+1\right)  -1},\ldots,f_{m+1}\right)  =C\left(
f_{m+1}\right)  =f_{m+1}%
\]
(by (\ref{eq.exa.ind.gen-ass-maps.Cex.1}), applied to $X_{m+1}$, $X_{m+2}$ and
$f_{m+1}$ instead of $X_{1}$, $X_{2}$ and $f_{1}$), so that%
\[
\underbrace{C\left(  f_{m+1},f_{\left(  m+1\right)  -1},\ldots,f_{k+1}\right)
}_{=f_{m+1}}\circ\underbrace{C\left(  f_{k},f_{k-1},\ldots,f_{1}\right)
}_{\substack{=C\left(  f_{m},f_{m-1},\ldots,f_{1}\right)  \\\text{(since
}k=m\text{)}}}=f_{m+1}\circ C\left(  f_{m},f_{m-1},\ldots,f_{1}\right)  .
\]
Comparing this with (\ref{pf.prop.ind.gen-ass-maps.Ceq.c1.pf.1}), we obtain%
\[
C\left(  f_{m+1},f_{\left(  m+1\right)  -1},\ldots,f_{1}\right)  =C\left(
f_{m+1},f_{\left(  m+1\right)  -1},\ldots,f_{k+1}\right)  \circ C\left(
f_{k},f_{k-1},\ldots,f_{1}\right)  .
\]
Hence, Claim 1 is proven in Case 1.

Let us now consider Case 2. In this case, we have $k\neq m$. Combining
$k\in\left\{  1,2,\ldots,m\right\}  $ with $k\neq m$, we obtain%
\[
k\in\left\{  1,2,\ldots,m\right\}  \setminus\left\{  m\right\}  =\left\{
1,2,\ldots,m-1\right\}  .
\]
Hence, $k\leq m-1<m$, so that $m+1-\underbrace{k}_{<m}>m+1-m=1$.

But we assumed that Proposition \ref{prop.ind.gen-ass-maps.Ceq} holds under
the condition that $n=m$. Hence, we can apply Proposition
\ref{prop.ind.gen-ass-maps.Ceq} to $m$ instead of $n$. We thus obtain%
\[
C\left(  f_{m},f_{m-1},\ldots,f_{1}\right)  =C\left(  f_{m},f_{m-1}%
,\ldots,f_{k+1}\right)  \circ C\left(  f_{k},f_{k-1},\ldots,f_{1}\right)
\]
(since $k\in\left\{  1,2,\ldots,m-1\right\}  $). Now,
(\ref{pf.prop.ind.gen-ass-maps.Ceq.c1.pf.1}) yields%
\begin{align}
&  C\left(  f_{m+1},f_{\left(  m+1\right)  -1},\ldots,f_{1}\right) \nonumber\\
&  =f_{m+1}\circ\underbrace{C\left(  f_{m},f_{m-1},\ldots,f_{1}\right)
}_{=C\left(  f_{m},f_{m-1},\ldots,f_{k+1}\right)  \circ C\left(  f_{k}%
,f_{k-1},\ldots,f_{1}\right)  }\nonumber\\
&  =f_{m+1}\circ\left(  C\left(  f_{m},f_{m-1},\ldots,f_{k+1}\right)  \circ
C\left(  f_{k},f_{k-1},\ldots,f_{1}\right)  \right)  .
\label{pf.prop.ind.gen-ass-maps.Ceq.c1.pf.4}%
\end{align}

On the other hand, $m+1-k>1$. Hence, (\ref{eq.def.ind.gen-ass-maps.C.rec})
(applied to $m+1-k$, $X_{k+i}$ and $f_{k+i}$ instead of $n$, $X_{i}$ and
$f_{i}$) yields%
\begin{align*}
&  C\left(  f_{k+\left(  m+1-k\right)  },f_{k+\left(  \left(  m+1-k\right)
-1\right)  },\ldots,f_{k+1}\right) \\
&  =f_{k+\left(  m+1-k\right)  }\circ C\left(  f_{k+\left(  \left(
m+1-k\right)  -1\right)  },f_{k+\left(  \left(  m+1-k\right)  -2\right)
},\ldots,f_{k+1}\right) \\
&  =f_{m+1}\circ C\left(  f_{m},f_{m-1},\ldots,f_{k+1}\right) \\
&  \ \ \ \ \ \ \ \ \ \ \left(
\begin{array}
[c]{c}%
\text{since }k+\left(  m+1-k\right)  =m+1\text{ and }k+\left(  \left(
m+1-k\right)  -1\right)  =m\\
\text{and }k+\left(  \left(  m+1-k\right)  -2\right)  =m-1
\end{array}
\right)  .
\end{align*}
Since $k+\left(  m+1-k\right)  =m+1$ and $k+\left(  \left(  m+1-k\right)
-1\right)  =\left(  m+1\right)  -1$, this rewrites as%
\[
C\left(  f_{m+1},f_{\left(  m+1\right)  -1},\ldots,f_{k+1}\right)
=f_{m+1}\circ C\left(  f_{m},f_{m-1},\ldots,f_{k+1}\right)  .
\]
Hence,%
\begin{align*}
&  \underbrace{C\left(  f_{m+1},f_{\left(  m+1\right)  -1},\ldots
,f_{k+1}\right)  }_{=f_{m+1}\circ C\left(  f_{m},f_{m-1},\ldots,f_{k+1}%
\right)  }\circ C\left(  f_{k},f_{k-1},\ldots,f_{1}\right) \\
&  =\left(  f_{m+1}\circ C\left(  f_{m},f_{m-1},\ldots,f_{k+1}\right)
\right)  \circ C\left(  f_{k},f_{k-1},\ldots,f_{1}\right) \\
&  =f_{m+1}\circ\left(  C\left(  f_{m},f_{m-1},\ldots,f_{k+1}\right)  \circ
C\left(  f_{k},f_{k-1},\ldots,f_{1}\right)  \right)
\end{align*}
(by Proposition \ref{prop.ind.gen-ass-maps.fgh}, applied to $X=X_{1}$,
$Y=X_{k+1}$, $Z=X_{m+1}$, $W=X_{m+2}$, $c=C\left(  f_{k},f_{k-1},\ldots
,f_{1}\right)  $, $b=C\left(  f_{m},f_{m-1},\ldots,f_{k+1}\right)  $ and
$a=f_{m+1}$). Comparing this with (\ref{pf.prop.ind.gen-ass-maps.Ceq.c1.pf.4}%
), we obtain%
\[
C\left(  f_{m+1},f_{\left(  m+1\right)  -1},\ldots,f_{1}\right)  =C\left(
f_{m+1},f_{\left(  m+1\right)  -1},\ldots,f_{k+1}\right)  \circ C\left(
f_{k},f_{k-1},\ldots,f_{1}\right)  .
\]
Hence, Claim 1 is proven in Case 2.

We have now proven Claim 1 in each of the two Cases 1 and 2. Since these two
Cases cover all possibilities, we thus conclude that Claim 1 always holds.]

Now, we have proven Claim 1. In other words, we have proven that Proposition
\ref{prop.ind.gen-ass-maps.Ceq} holds under the condition that $n=m+1$. This
completes the induction step. Hence, Proposition
\ref{prop.ind.gen-ass-maps.Ceq} is proven by induction.
\end{proof}

\subsubsection{Proof of general associativity}

\begin{proposition}
\label{prop.ind.gen-ass-maps.Ceq-cp}Let $n$ be a positive integer. Let
$X_{1},X_{2},\ldots,X_{n+1}$ be $n+1$ sets. For each $i\in\left\{
1,2,\ldots,n\right\}  $, let $f_{i}:X_{i}\rightarrow X_{i+1}$ be a map. Then,
every complete parenthesization of $f_{n}\circ f_{n-1}\circ\cdots\circ f_{1}$
equals $C\left(  f_{n},f_{n-1},\ldots,f_{1}\right)  $.
\end{proposition}

\begin{proof}
[Proof of Proposition \ref{prop.ind.gen-ass-maps.Ceq-cp}.]Forget that we fixed
$n$, $X_{1},X_{2},\ldots,X_{n+1}$ and the maps $f_{i}$. We shall prove
Proposition \ref{prop.ind.gen-ass-maps.Ceq-cp} by strong induction on
$n$:\ \ \ \ \footnote{The induction principle that we are applying here is
Theorem \ref{thm.ind.SIP} with $g=1$ (since $\mathbb{Z}_{\geq1}$ is the set of
all positive integers).}

\textit{Induction step:} Let $m\in\mathbb{Z}_{\geq1}$. Assume that Proposition
\ref{prop.ind.gen-ass-maps.Ceq-cp} holds under the condition that $n<m$. We
must prove that Proposition \ref{prop.ind.gen-ass-maps.Ceq-cp} holds under the
condition that $n=m$. In other words, we must prove the following claim:

\begin{statement}
\textit{Claim 1:} Let $X_{1},X_{2},\ldots,X_{m+1}$ be $m+1$ sets. For each
$i\in\left\{  1,2,\ldots,m\right\}  $, let $f_{i}:X_{i}\rightarrow X_{i+1}$ be
a map. Then, every complete parenthesization of $f_{m}\circ f_{m-1}\circ
\cdots\circ f_{1}$ equals $C\left(  f_{m},f_{m-1},\ldots,f_{1}\right)  $.
\end{statement}

[\textit{Proof of Claim 1:} Let $\gamma$ be a complete parenthesization of
$f_{m}\circ f_{m-1}\circ\cdots\circ f_{1}$. Thus, we must prove that
$\gamma=C\left(  f_{m},f_{m-1},\ldots,f_{1}\right)  $.

We have $m\in\mathbb{Z}_{\geq1}$, thus $m\geq1$. Hence, either $m=1$ or $m>1$.
Thus, we are in one of the following two cases:

\textit{Case 1:} We have $m=1$.

\textit{Case 2:} We have $m>1$.

Let us first consider Case 1. In this case, we have $m=1$. Thus, we have
$C\left(  f_{m},f_{m-1},\ldots,f_{1}\right)  =f_{1}$ (by the definition of
$C\left(  f_{m},f_{m-1},\ldots,f_{1}\right)  $).

Recall that $m=1$. Thus, the definition of a \textquotedblleft complete
parenthesization of $f_{m}\circ f_{m-1}\circ\cdots\circ f_{1}$%
\textquotedblright\ shows that there is only one complete parenthesization of
$f_{m}\circ f_{m-1}\circ\cdots\circ f_{1}$, and this is simply the map
$f_{1}:X_{1}\rightarrow X_{2}$. Hence, $\gamma$ is simply the map $f_{1}%
:X_{1}\rightarrow X_{2}$ (since $\gamma$ is a complete parenthesization of
$f_{m}\circ f_{m-1}\circ\cdots\circ f_{1}$). Thus, $\gamma=f_{1}=C\left(
f_{m},f_{m-1},\ldots,f_{1}\right)  $ (since $C\left(  f_{m},f_{m-1}%
,\ldots,f_{1}\right)  =f_{1}$). Thus, $\gamma=C\left(  f_{m},f_{m-1}%
,\ldots,f_{1}\right)  $ is proven in Case 1.

Now, let us consider Case 2. In this case, we have $m>1$. Hence, the
definition of a \textquotedblleft complete parenthesization of $f_{m}\circ
f_{m-1}\circ\cdots\circ f_{1}$\textquotedblright\ shows that any complete
parenthesization of $f_{m}\circ f_{m-1}\circ\cdots\circ f_{1}$ is a map of the
form $\alpha\circ\beta$, where

\begin{itemize}
\item $k$ is some element of $\left\{  1,2,\ldots,m-1\right\}  $;

\item $\alpha$ is a complete parenthesization of $f_{m}\circ f_{m-1}%
\circ\cdots\circ f_{k+1}$;

\item $\beta$ is a complete parenthesization of $f_{k}\circ f_{k-1}\circ
\cdots\circ f_{1}$.
\end{itemize}

Thus, $\gamma$ is a map of this form (since $\gamma$ is a complete
parenthesization of $f_{m}\circ f_{m-1}\circ\cdots\circ f_{1}$). In other
words, we can write $\gamma$ in the form $\gamma=\alpha\circ\beta$, where $k$
is some element of $\left\{  1,2,\ldots,m-1\right\}  $, where $\alpha$ is a
complete parenthesization of $f_{m}\circ f_{m-1}\circ\cdots\circ f_{k+1}$, and
where $\beta$ is a complete parenthesization of $f_{k}\circ f_{k-1}\circ
\cdots\circ f_{1}$. Consider these $k$, $\alpha$ and $\beta$.

We have $k\in\left\{  1,2,\ldots,m-1\right\}  $, thus $k\leq m-1<m$. Hence, we
can apply Proposition \ref{prop.ind.gen-ass-maps.Ceq-cp} to $n=k$ (since we
assumed that Proposition \ref{prop.ind.gen-ass-maps.Ceq-cp} holds under the
condition that $n<m$). We thus conclude that every complete parenthesization
of $f_{k}\circ f_{k-1}\circ\cdots\circ f_{1}$ equals $C\left(  f_{k}%
,f_{k-1},\ldots,f_{1}\right)  $. Hence, $\beta$ equals $C\left(  f_{k}%
,f_{k-1},\ldots,f_{1}\right)  $ (since $\beta$ is a complete parenthesization
of $f_{k}\circ f_{k-1}\circ\cdots\circ f_{1}$). In other words,%
\begin{equation}
\beta=C\left(  f_{k},f_{k-1},\ldots,f_{1}\right)  .
\label{pf.prop.ind.gen-ass-maps.Ceq-cp.4b}%
\end{equation}

We have $k\in\left\{  1,2,\ldots,m-1\right\}  $, thus $k\geq1$ and therefore
$m-\underbrace{k}_{\geq1}\leq m-1<m$. Hence, we can apply Proposition
\ref{prop.ind.gen-ass-maps.Ceq-cp} to $m-k$, $X_{k+i}$ and $f_{k+i}$ instead
of $n$, $X_{i}$ and $f_{i}$ (since we assumed that Proposition
\ref{prop.ind.gen-ass-maps.Ceq-cp} holds under the condition that $n<m$). We
thus conclude that every complete parenthesization of $f_{k+\left(
m-k\right)  }\circ f_{k+\left(  m-k-1\right)  }\circ\cdots\circ f_{k+1}$
equals $C\left(  f_{k+\left(  m-k\right)  },f_{k+\left(  m-k-1\right)
},\ldots,f_{k+1}\right)  $.

Since $k+\left(  m-k\right)  =m$ and $k+\left(  m-k-1\right)  =m-1$, this
rewrites as follows: Every complete parenthesization of $f_{m}\circ
f_{m-1}\circ\cdots\circ f_{k+1}$ equals $C\left(  f_{m},f_{m-1},\ldots
,f_{k+1}\right)  $. Thus, $\alpha$ equals $C\left(  f_{m},f_{m-1}%
,\ldots,f_{k+1}\right)  $ (since $\alpha$ is a complete parenthesization of
$f_{m}\circ f_{m-1}\circ\cdots\circ f_{k+1}$). In other words,%
\begin{equation}
\alpha=C\left(  f_{m},f_{m-1},\ldots,f_{k+1}\right)  .
\label{pf.prop.ind.gen-ass-maps.Ceq-cp.4a}%
\end{equation}

But Proposition \ref{prop.ind.gen-ass-maps.Ceq} (applied to $n=m$) yields%
\[
C\left(  f_{m},f_{m-1},\ldots,f_{1}\right)  =\underbrace{C\left(
f_{m},f_{m-1},\ldots,f_{k+1}\right)  }_{\substack{=\alpha\\\text{(by
(\ref{pf.prop.ind.gen-ass-maps.Ceq-cp.4a}))}}}\circ\underbrace{C\left(
f_{k},f_{k-1},\ldots,f_{1}\right)  }_{\substack{=\beta\\\text{(by
(\ref{pf.prop.ind.gen-ass-maps.Ceq-cp.4b}))}}}=\alpha\circ\beta=\gamma
\]
(since $\gamma=\alpha\circ\beta$), so that $\gamma=C\left(  f_{m}%
,f_{m-1},\ldots,f_{1}\right)  $. Hence, $\gamma=C\left(  f_{m},f_{m-1}%
,\ldots,f_{1}\right)  $ is proven in Case 2.

We now have shown that $\gamma=C\left(  f_{m},f_{m-1},\ldots,f_{1}\right)  $
in each of the two Cases 1 and 2. Since these two Cases cover all
possibilities, this yields that $\gamma=C\left(  f_{m},f_{m-1},\ldots
,f_{1}\right)  $ always holds.

Now, forget that we fixed $\gamma$. We thus have shown that $\gamma=C\left(
f_{m},f_{m-1},\ldots,f_{1}\right)  $ whenever $\gamma$ is a complete
parenthesization of $f_{m}\circ f_{m-1}\circ\cdots\circ f_{1}$. In other
words, every complete parenthesization of $f_{m}\circ f_{m-1}\circ\cdots\circ
f_{1}$ equals $C\left(  f_{m},f_{m-1},\ldots,f_{1}\right)  $. This proves
Claim 1.]

Now, we have proven Claim 1. In other words, we have proven that Proposition
\ref{prop.ind.gen-ass-maps.Ceq-cp} holds under the condition that $n=m$. This
completes the induction step. Hence, Proposition
\ref{prop.ind.gen-ass-maps.Ceq-cp} is proven by strong induction.
\end{proof}

\begin{proof}
[Proof of Theorem \ref{thm.ind.gen-ass-maps.cp}.]Proposition
\ref{prop.ind.gen-ass-maps.Ceq-cp} shows that every complete parenthesization
of $f_{n}\circ f_{n-1}\circ\cdots\circ f_{1}$ equals $C\left(  f_{n}%
,f_{n-1},\ldots,f_{1}\right)  $. Thus, all complete parenthesizations of
$f_{n}\circ f_{n-1}\circ\cdots\circ f_{1}$ are the same map. This proves
Theorem \ref{thm.ind.gen-ass-maps.cp}.
\end{proof}

\subsubsection{Compositions of multiple maps without parentheses}

\begin{definition}
\label{def.ind.gen-ass-maps.comp}Let $n$ be a positive integer. Let
$X_{1},X_{2},\ldots,X_{n+1}$ be $n+1$ sets. For each $i\in\left\{
1,2,\ldots,n\right\}  $, let $f_{i}:X_{i}\rightarrow X_{i+1}$ be a map. Then,
the map $C\left(  f_{n},f_{n-1},\ldots,f_{1}\right)  :X_{1}\rightarrow
X_{n+1}$ is denoted by $f_{n}\circ f_{n-1}\circ\cdots\circ f_{1}$ and called
the \textit{composition} of $f_{n},f_{n-1},\ldots,f_{1}$. This notation
$f_{n}\circ f_{n-1}\circ\cdots\circ f_{1}$ may conflict with existing
notations in two cases:

\begin{itemize}
\item In the case when $n=1$, this notation $f_{n}\circ f_{n-1}\circ
\cdots\circ f_{1}$ simply becomes $f_{1}$, which looks exactly like the map
$f_{1}$ itself. Fortunately, this conflict of notation is harmless, because
the new meaning that we are giving to $f_{1}$ in this case (namely, $C\left(
f_{n},f_{n-1},\ldots,f_{1}\right)  =C\left(  f_{1}\right)  $) agrees with the
map $f_{1}$ (because of (\ref{eq.exa.ind.gen-ass-maps.Cex.1})).

\item In the case when $n=2$, this notation $f_{n}\circ f_{n-1}\circ
\cdots\circ f_{1}$ simply becomes $f_{2}\circ f_{1}$, which looks exactly like
the composition $f_{2}\circ f_{1}$ of the two maps $f_{2}$ and $f_{1}$.
Fortunately, this conflict of notation is harmless, because the new meaning
that we are giving to $f_{2}\circ f_{1}$ in this case (namely, $C\left(
f_{n},f_{n-1},\ldots,f_{1}\right)  =C\left(  f_{2},f_{1}\right)  $) agrees
with the latter composition (because of (\ref{eq.exa.ind.gen-ass-maps.Cex.2})).
\end{itemize}

Thus, in both cases, the conflict with existing notations is harmless (the
conflicting notations actually stand for the same thing).
\end{definition}

\begin{remark}
\label{rmk.ind.gen-ass-maps.drop}Let $n$, $X_{1},X_{2},\ldots,X_{n+1}$ and
$f_{i}$ be as in Definition \ref{def.ind.gen-ass-maps.comp}. Then, Proposition
\ref{prop.ind.gen-ass-maps.Ceq-cp} shows that every complete parenthesization
of $f_{n}\circ f_{n-1}\circ\cdots\circ f_{1}$ equals $C\left(  f_{n}%
,f_{n-1},\ldots,f_{1}\right)  $. In other words, every complete
parenthesization of $f_{n}\circ f_{n-1}\circ\cdots\circ f_{1}$ equals
$f_{n}\circ f_{n-1}\circ\cdots\circ f_{1}$ (because $f_{n}\circ f_{n-1}%
\circ\cdots\circ f_{1}$ was defined to be $C\left(  f_{n},f_{n-1},\ldots
,f_{1}\right)  $ in Definition \ref{def.ind.gen-ass-maps.comp}). In other
words, \textbf{we can drop the parentheses} in every complete parenthesization
of $f_{n}\circ f_{n-1}\circ\cdots\circ f_{1}$. For example, for $n=7$, we get%
\[
\left(  f_{7}\circ\left(  f_{6}\circ f_{5}\right)  \right)  \circ\left(
f_{4}\circ\left(  \left(  f_{3}\circ f_{2}\right)  \circ f_{1}\right)
\right)  =f_{7}\circ f_{6}\circ f_{5}\circ f_{4}\circ f_{3}\circ f_{2}\circ
f_{1},
\]
and a similar equality for any other complete parenthesization.
\end{remark}

Definition \ref{def.ind.gen-ass-maps.comp} and Remark
\ref{rmk.ind.gen-ass-maps.drop} finally give us the justification to write
compositions of multiple maps (like $f_{n}\circ f_{n-1}\circ\cdots\circ f_{1}$
for $n\geq1$) without the need for parentheses. We shall now go one little
step further and extend this notation to the case of $n=0$ -- that is, we
shall define the composition of \textbf{no} maps:

\begin{definition}
\label{def.ind.gen-ass-maps.comp0}Let $n\in\mathbb{N}$. Let $X_{1}%
,X_{2},\ldots,X_{n+1}$ be $n+1$ sets. For each $i\in\left\{  1,2,\ldots
,n\right\}  $, let $f_{i}:X_{i}\rightarrow X_{i+1}$ be a map. In Definition
\ref{def.ind.gen-ass-maps.comp}, we have defined the composition $f_{n}\circ
f_{n-1}\circ\cdots\circ f_{1}$ of $f_{n},f_{n-1},\ldots,f_{1}$ when $n$ is a
positive integer. We shall now extend this definition to the case when $n=0$
(so that it will be defined for all $n\in\mathbb{N}$, not just for all
positive integers $n$). Namely, we extend it by setting
\begin{equation}
f_{n}\circ f_{n-1}\circ\cdots\circ f_{1}=\operatorname*{id}\nolimits_{X_{1}%
}\ \ \ \ \ \ \ \ \ \ \text{when }n=0. \label{eq.def.ind.gen-ass-maps.comp0.0}%
\end{equation}
That is, we say that the composition of $0$ maps is the identity map
$\operatorname*{id}\nolimits_{X_{1}}:X_{1}\rightarrow X_{1}$. This composition
of $0$ maps is also known as the \textit{empty composition of maps}. Thus, the
empty composition of maps is defined to be $\operatorname*{id}\nolimits_{X_{1}%
}$. (This is similar to the well-known conventions that a sum of $0$ numbers
is $0$, and that a product of $0$ numbers is $1$.)

This definition is slightly dangerous, because it entails that the composition
of $0$ maps depends on the set $X_{1}$, but of course the $0$ maps being
composed know nothing about this set $X_{1}$. Thus, when we speak of an empty
composition, we should always specify the set $X_{1}$ or ensure that it is
clear from the context. (See Definition \ref{def.ind.gen-ass-maps.fn} below
for an example where it is clear from the context.)
\end{definition}

\subsubsection{Composition powers}

Having defined the composition of $n$ maps, we get the notion of composition
powers of maps for free:

\begin{definition}
\label{def.ind.gen-ass-maps.fn}Let $n\in\mathbb{N}$. Let $X$ be a set. Let
$f:X\rightarrow X$ be a map. Then, $f^{\circ n}$ shall denote the map%
\[
\underbrace{f\circ f\circ\cdots\circ f}_{n\text{ times }f}:X\rightarrow X.
\]
Thus, in particular,
\begin{equation}
f^{\circ0}=\underbrace{f\circ f\circ\cdots\circ f}_{0\text{ times }%
f}=\operatorname*{id}\nolimits_{X} \label{eq.def.ind.gen-ass-maps.fn.f0}%
\end{equation}
(by Definition \ref{def.ind.gen-ass-maps.comp0}). Also, $f^{\circ1}=f$,
$f^{\circ2}=f\circ f$, $f^{\circ3}=f\circ f\circ f$, etc.

The map $f^{\circ n}$ is called the $n$\textit{-th composition power} of $f$
(or simply the $n$\textit{-th power} of $f$).
\end{definition}

Before we study composition powers in detail, let us show a general rule that
allows us to \textquotedblleft split\textquotedblright\ compositions of maps:

\begin{theorem}
\label{thm.ind.gen-ass-maps.1}Let $n\in\mathbb{N}$. Let $X_{1},X_{2}%
,\ldots,X_{n+1}$ be $n+1$ sets. For each $i\in\left\{  1,2,\ldots,n\right\}
$, let $f_{i}:X_{i}\rightarrow X_{i+1}$ be a map.

\textbf{(a)} We have
\[
f_{n}\circ f_{n-1}\circ\cdots\circ f_{1}=\left(  f_{n}\circ f_{n-1}\circ
\cdots\circ f_{k+1}\right)  \circ\left(  f_{k}\circ f_{k-1}\circ\cdots\circ
f_{1}\right)
\]
for each $k\in\left\{  0,1,\ldots,n\right\}  $.

\textbf{(b)} If $n\geq1$, then%
\[
f_{n}\circ f_{n-1}\circ\cdots\circ f_{1}=f_{n}\circ\left(  f_{n-1}\circ
f_{n-2}\circ\cdots\circ f_{1}\right)  .
\]

\textbf{(c)} If $n\geq1$, then%
\[
f_{n}\circ f_{n-1}\circ\cdots\circ f_{1}=\left(  f_{n}\circ f_{n-1}\circ
\cdots\circ f_{2}\right)  \circ f_{1}.
\]

\end{theorem}

\begin{proof}
[Proof of Theorem \ref{thm.ind.gen-ass-maps.1}.]\textbf{(a)} Let $k\in\left\{
0,1,\ldots,n\right\}  $. We are in one of the following three cases:

\textit{Case 1:} We have $k=0$.

\textit{Case 2:} We have $k=n$.

\textit{Case 3:} We have neither $k=0$ nor $k=n$.

(Of course, Cases 1 and 2 overlap when $n=0$.)

Let us first consider Case 1. In this case, we have $k=0$. Thus, $f_{k}\circ
f_{k-1}\circ\cdots\circ f_{1}$ is an empty composition of maps, and therefore
equals $\operatorname*{id}\nolimits_{X_{1}}$. In other words, $f_{k}\circ
f_{k-1}\circ\cdots\circ f_{1}=\operatorname*{id}\nolimits_{X_{1}}$.

On the other hand, $k=0$, so that $k+1=1$. Hence, $f_{n}\circ f_{n-1}%
\circ\cdots\circ f_{k+1}=f_{n}\circ f_{n-1}\circ\cdots\circ f_{1}$. Thus,%
\begin{align*}
&  \underbrace{\left(  f_{n}\circ f_{n-1}\circ\cdots\circ f_{k+1}\right)
}_{=f_{n}\circ f_{n-1}\circ\cdots\circ f_{1}}\circ\underbrace{\left(
f_{k}\circ f_{k-1}\circ\cdots\circ f_{1}\right)  }_{=\operatorname*{id}%
\nolimits_{X_{1}}}\\
&  =\left(  f_{n}\circ f_{n-1}\circ\cdots\circ f_{1}\right)  \circ
\operatorname*{id}\nolimits_{X_{1}}=f_{n}\circ f_{n-1}\circ\cdots\circ f_{1}.
\end{align*}
In other words,%
\[
f_{n}\circ f_{n-1}\circ\cdots\circ f_{1}=\left(  f_{n}\circ f_{n-1}\circ
\cdots\circ f_{k+1}\right)  \circ\left(  f_{k}\circ f_{k-1}\circ\cdots\circ
f_{1}\right)  .
\]
Hence, Theorem \ref{thm.ind.gen-ass-maps.1} \textbf{(a)} is proven in Case 1.

Let us next consider Case 2. In this case, we have $k=n$. Thus, $f_{n}\circ
f_{n-1}\circ\cdots\circ f_{k+1}$ is an empty composition of maps, and
therefore equals $\operatorname*{id}\nolimits_{X_{n+1}}$. In other words,
$f_{n}\circ f_{n-1}\circ\cdots\circ f_{k+1}=\operatorname*{id}%
\nolimits_{X_{n+1}}$. Thus,%
\begin{align*}
&  \underbrace{\left(  f_{n}\circ f_{n-1}\circ\cdots\circ f_{k+1}\right)
}_{=\operatorname*{id}\nolimits_{X_{n+1}}}\circ\underbrace{\left(  f_{k}\circ
f_{k-1}\circ\cdots\circ f_{1}\right)  }_{\substack{=f_{n}\circ f_{n-1}%
\circ\cdots\circ f_{1}\\\text{(since }k=n\text{)}}}\\
&  =\operatorname*{id}\nolimits_{X_{n+1}}\circ\left(  f_{n}\circ f_{n-1}%
\circ\cdots\circ f_{1}\right)  =f_{n}\circ f_{n-1}\circ\cdots\circ f_{1}.
\end{align*}
In other words,%
\[
f_{n}\circ f_{n-1}\circ\cdots\circ f_{1}=\left(  f_{n}\circ f_{n-1}\circ
\cdots\circ f_{k+1}\right)  \circ\left(  f_{k}\circ f_{k-1}\circ\cdots\circ
f_{1}\right)  .
\]
Hence, Theorem \ref{thm.ind.gen-ass-maps.1} \textbf{(a)} is proven in Case 2.

Let us finally consider Case 3. In this case, we have neither $k=0$ nor $k=n$.
In other words, we have $k\neq0$ and $k\neq n$. Combining $k\in\left\{
0,1,\ldots,n\right\}  $ with $k\neq0$, we find $k\in\left\{  0,1,\ldots
,n\right\}  \setminus\left\{  0\right\}  \subseteq\left\{  1,2,\ldots
,n\right\}  $. Combining this with $k\neq n$, we find $k\in\left\{
1,2,\ldots,n\right\}  \setminus\left\{  n\right\}  \subseteq\left\{
1,2,\ldots,n-1\right\}  $. Hence, $1\leq k\leq n-1$, so that $n-1\geq1$ and
thus $n\geq2\geq1$. Hence, $n$ is a positive integer. Thus, Proposition
\ref{prop.ind.gen-ass-maps.Ceq} yields
\begin{equation}
C\left(  f_{n},f_{n-1},\ldots,f_{1}\right)  =C\left(  f_{n},f_{n-1}%
,\ldots,f_{k+1}\right)  \circ C\left(  f_{k},f_{k-1},\ldots,f_{1}\right)  .
\label{pf.thm.ind.gen-ass-maps.1.a.2}%
\end{equation}

Now, $k$ is a positive integer (since $k\in\left\{  1,2,\ldots,n-1\right\}
$). Hence, Definition \ref{def.ind.gen-ass-maps.comp} (applied to $k$ instead
of $n$) yields%
\begin{equation}
f_{k}\circ f_{k-1}\circ\cdots\circ f_{1}=C\left(  f_{k},f_{k-1},\ldots
,f_{1}\right)  . \label{pf.thm.ind.gen-ass-maps.1.a.3a}%
\end{equation}
Also, $n-k$ is an integer satisfying $n-\underbrace{k}_{\leq n-1}\geq
n-\left(  n-1\right)  =1$. Hence, $n-k$ is a positive integer. Thus,
Definition \ref{def.ind.gen-ass-maps.comp} (applied to $n-k$, $X_{k+i}$ and
$f_{k+i}$ instead of $n$, $X_{i}$ and $f_{i}$) yields%
\[
f_{k+\left(  n-k\right)  }\circ f_{k+\left(  n-k-1\right)  }\circ\cdots\circ
f_{k+1}=C\left(  f_{k+\left(  n-k\right)  },f_{k+\left(  n-k-1\right)
},\ldots,f_{k+1}\right)  .
\]
In view of $k+\left(  n-k\right)  =n$ and $k+\left(  n-k-1\right)  =n-1$, this
rewrites as follows:%
\begin{equation}
f_{n}\circ f_{n-1}\circ\cdots\circ f_{k+1}=C\left(  f_{n},f_{n-1}%
,\ldots,f_{k+1}\right)  . \label{pf.thm.ind.gen-ass-maps.1.a.3b}%
\end{equation}

But $n$ is a positive integer. Thus, Definition
\ref{def.ind.gen-ass-maps.comp} yields
\begin{align*}
f_{n}\circ f_{n-1}\circ\cdots\circ f_{1}  &  =C\left(  f_{n},f_{n-1}%
,\ldots,f_{1}\right) \\
&  =\underbrace{C\left(  f_{n},f_{n-1},\ldots,f_{k+1}\right)  }%
_{\substack{=f_{n}\circ f_{n-1}\circ\cdots\circ f_{k+1}\\\text{(by
(\ref{pf.thm.ind.gen-ass-maps.1.a.3b}))}}}\circ\underbrace{C\left(
f_{k},f_{k-1},\ldots,f_{1}\right)  }_{\substack{=f_{k}\circ f_{k-1}\circ
\cdots\circ f_{1}\\\text{(by (\ref{pf.thm.ind.gen-ass-maps.1.a.3a}))}%
}}\ \ \ \ \ \ \ \ \ \ \left(  \text{by (\ref{pf.thm.ind.gen-ass-maps.1.a.2}%
)}\right) \\
&  =\left(  f_{n}\circ f_{n-1}\circ\cdots\circ f_{k+1}\right)  \circ\left(
f_{k}\circ f_{k-1}\circ\cdots\circ f_{1}\right)  .
\end{align*}
Hence, Theorem \ref{thm.ind.gen-ass-maps.1} \textbf{(a)} is proven in Case 3.

We have now proven Theorem \ref{thm.ind.gen-ass-maps.1} \textbf{(a)} in each
of the three Cases 1, 2 and 3. Since these three Cases cover all
possibilities, we thus conclude that Theorem \ref{thm.ind.gen-ass-maps.1}
\textbf{(a)} always holds.

\textbf{(b)} Assume that $n\geq1$. Hence, $n-1\geq0$, so that $n-1\in\left\{
0,1,\ldots,n\right\}  $ (since $n-1\leq n$). Hence, Theorem
\ref{thm.ind.gen-ass-maps.1} \textbf{(a)} (applied to $k=n-1$) yields%
\begin{align*}
f_{n}\circ f_{n-1}\circ\cdots\circ f_{1}  &  =\underbrace{\left(  f_{n}\circ
f_{n-1}\circ\cdots\circ f_{\left(  n-1\right)  +1}\right)  }_{=f_{n}\circ
f_{n-1}\circ\cdots\circ f_{n}=f_{n}}\circ\underbrace{\left(  f_{n-1}\circ
f_{\left(  n-1\right)  -1}\circ\cdots\circ f_{1}\right)  }_{=f_{n-1}\circ
f_{n-2}\circ\cdots\circ f_{1}}\\
&  =f_{n}\circ\left(  f_{n-1}\circ f_{n-2}\circ\cdots\circ f_{1}\right)  .
\end{align*}
This proves Theorem \ref{thm.ind.gen-ass-maps.1} \textbf{(b)}.

\textbf{(c)} Assume that $n\geq1$. Hence, $1\in\left\{  0,1,\ldots,n\right\}
$. Hence, Theorem \ref{thm.ind.gen-ass-maps.1} \textbf{(a)} (applied to $k=1$)
yields%
\begin{align*}
f_{n}\circ f_{n-1}\circ\cdots\circ f_{1}  &  =\underbrace{\left(  f_{n}\circ
f_{n-1}\circ\cdots\circ f_{1+1}\right)  }_{=f_{n}\circ f_{n-1}\circ\cdots\circ
f_{2}}\circ\underbrace{\left(  f_{1}\circ f_{1-1}\circ\cdots\circ
f_{1}\right)  }_{=f_{1}}\\
&  =\left(  f_{n}\circ f_{n-1}\circ\cdots\circ f_{2}\right)  \circ f_{1}.
\end{align*}
This proves Theorem \ref{thm.ind.gen-ass-maps.1} \textbf{(c)}.
\end{proof}

We can draw some consequences about composition powers of maps from Theorem
\ref{thm.ind.gen-ass-maps.1}:

\begin{proposition}
\label{prop.ind.map-powers.1}Let $X$ be a set. Let $f:X\rightarrow X$ be a
map. Let $n$ be a positive integer.

\textbf{(a)} We have $f^{\circ n}=f\circ f^{\circ\left(  n-1\right)  }$.

\textbf{(b)} We have $f^{\circ n}=f^{\circ\left(  n-1\right)  }\circ f$.
\end{proposition}

\begin{proof}
[Proof of Proposition \ref{prop.ind.map-powers.1}.]The definition of $f^{\circ
n}$ yields%
\begin{equation}
f^{\circ n}=\underbrace{f\circ f\circ\cdots\circ f}_{n\text{ times }f}.
\label{pf.prop.ind.map-powers.1.a.1}%
\end{equation}
The definition of $f^{\circ\left(  n-1\right)  }$ yields%
\begin{equation}
f^{\circ\left(  n-1\right)  }=\underbrace{f\circ f\circ\cdots\circ
f}_{n-1\text{ times }f}. \label{pf.prop.ind.map-powers.1.a.2}%
\end{equation}

\textbf{(a)} Theorem \ref{thm.ind.gen-ass-maps.1} \textbf{(b)} (applied to
$X_{i}=X$ and $f_{i}=f$) yields%
\[
\underbrace{f\circ f\circ\cdots\circ f}_{n\text{ times }f}=f\circ\left(
\underbrace{f\circ f\circ\cdots\circ f}_{n-1\text{ times }f}\right)  .
\]
In view of (\ref{pf.prop.ind.map-powers.1.a.1}) and
(\ref{pf.prop.ind.map-powers.1.a.2}), this rewrites as $f^{\circ n}=f\circ
f^{\circ\left(  n-1\right)  }$. This proves Proposition
\ref{prop.ind.map-powers.1} \textbf{(a)}.

\textbf{(b)} Theorem \ref{thm.ind.gen-ass-maps.1} \textbf{(c)} (applied to
$X_{i}=X$ and $f_{i}=f$) yields%
\[
\underbrace{f\circ f\circ\cdots\circ f}_{n\text{ times }f}=\left(
\underbrace{f\circ f\circ\cdots\circ f}_{n-1\text{ times }f}\right)  \circ f.
\]
In view of (\ref{pf.prop.ind.map-powers.1.a.1}) and
(\ref{pf.prop.ind.map-powers.1.a.2}), this rewrites as $f^{\circ n}%
=f^{\circ\left(  n-1\right)  }\circ f$. This proves Proposition
\ref{prop.ind.map-powers.1} \textbf{(b)}.
\end{proof}

\begin{proposition}
\label{prop.ind.map-powers.2}Let $X$ be a set. Let $f:X\rightarrow X$ be a map.

\textbf{(a)} We have $f^{\circ\left(  a+b\right)  }=f^{\circ a}\circ f^{\circ
b}$ for every $a\in\mathbb{N}$ and $b\in\mathbb{N}$.

\textbf{(b)} We have $f^{\circ\left(  ab\right)  }=\left(  f^{\circ a}\right)
^{\circ b}$ for every $a\in\mathbb{N}$ and $b\in\mathbb{N}$.
\end{proposition}

\begin{proof}
[Proof of Proposition \ref{prop.ind.map-powers.2}.]\textbf{(a)} Let
$a\in\mathbb{N}$ and $b\in\mathbb{N}$. Thus, $a\geq0$ and $b\geq0$, so that
$0\leq b\leq a+b$ (since $\underbrace{a}_{\geq0}+b\geq b$). Hence,
$b\in\left\{  0,1,\ldots,a+b\right\}  $. Thus, Theorem
\ref{thm.ind.gen-ass-maps.1} \textbf{(a)} (applied to $n=a+b$, $X_{i}=X$,
$f_{i}=f$ and $k=b$) yields%
\[
\underbrace{f\circ f\circ\cdots\circ f}_{a+b\text{ times }f}=\left(
\underbrace{f\circ f\circ\cdots\circ f}_{\left(  a+b\right)  -b\text{ times
}f}\right)  \circ\left(  \underbrace{f\circ f\circ\cdots\circ f}_{b\text{
times }f}\right)  .
\]
In view of $\left(  a+b\right)  -b=a$, this rewrites as%
\begin{equation}
\underbrace{f\circ f\circ\cdots\circ f}_{a+b\text{ times }f}=\left(
\underbrace{f\circ f\circ\cdots\circ f}_{a\text{ times }f}\right)
\circ\left(  \underbrace{f\circ f\circ\cdots\circ f}_{b\text{ times }%
f}\right)  . \label{pf.prop.ind.map-powers.2.1}%
\end{equation}

But the definition of $f^{\circ a}$ yields%
\begin{equation}
f^{\circ a}=\underbrace{f\circ f\circ\cdots\circ f}_{a\text{ times }f}.
\label{pf.prop.ind.map-powers.2.2}%
\end{equation}
Also, the definition of $f^{\circ b}$ yields%
\begin{equation}
f^{\circ b}=\underbrace{f\circ f\circ\cdots\circ f}_{b\text{ times }f}.
\label{pf.prop.ind.map-powers.2.3}%
\end{equation}
Finally, the definition of $f^{\circ\left(  a+b\right)  }$ yields%
\begin{equation}
f^{\circ\left(  a+b\right)  }=\underbrace{f\circ f\circ\cdots\circ
f}_{a+b\text{ times }f}. \label{pf.prop.ind.map-powers.2.4}%
\end{equation}
In view of these three equalities (\ref{pf.prop.ind.map-powers.2.2}),
(\ref{pf.prop.ind.map-powers.2.3}) and (\ref{pf.prop.ind.map-powers.2.4}), we
can rewrite the equality (\ref{pf.prop.ind.map-powers.2.1}) as $f^{\circ
\left(  a+b\right)  }=f^{\circ a}\circ f^{\circ b}$. This proves Proposition
\ref{prop.ind.map-powers.2} \textbf{(a)}.

(Alternatively, it is easy to prove Proposition \ref{prop.ind.map-powers.2}
\textbf{(a)} by induction on $a$.)

\textbf{(b)} Let $a\in\mathbb{N}$. We claim that%
\begin{equation}
f^{\circ\left(  ab\right)  }=\left(  f^{\circ a}\right)  ^{\circ
b}\ \ \ \ \ \ \ \ \ \ \text{for every }b\in\mathbb{N}.
\label{pf.prop.ind.map-powers.2.goal}%
\end{equation}
We shall prove (\ref{pf.prop.ind.map-powers.2.goal}) by induction on $b$:

\textit{Induction base:} We have $a\cdot0=0$ and thus $f^{\circ\left(
a\cdot0\right)  }=f^{\circ0}=\operatorname*{id}\nolimits_{X}$ (by
(\ref{eq.def.ind.gen-ass-maps.fn.f0})). Comparing this with $\left(  f^{\circ
a}\right)  ^{\circ0}=\operatorname*{id}\nolimits_{X}$ (which follows from
(\ref{eq.def.ind.gen-ass-maps.fn.f0}), applied to $f^{\circ a}$ instead of
$f$), we obtain $f^{\circ\left(  a\cdot0\right)  }=\left(  f^{\circ a}\right)
^{\circ0}$. In other words, (\ref{pf.prop.ind.map-powers.2.goal}) holds for
$b=0$. This completes the induction base.

\textit{Induction step:} Let $m\in\mathbb{N}$. Assume that
(\ref{pf.prop.ind.map-powers.2.goal}) holds for $b=m$. We must prove that
(\ref{pf.prop.ind.map-powers.2.goal}) holds for $b=m+1$.

We have assumed that (\ref{pf.prop.ind.map-powers.2.goal}) holds for $b=m$. In
other words, we have $f^{\circ\left(  am\right)  }=\left(  f^{\circ a}\right)
^{\circ m}$.

But $m+1$ is a positive integer (since $m+1>m\geq0$). Hence, Proposition
\ref{prop.ind.map-powers.1} \textbf{(b)} (applied to $m+1$ and $f^{\circ a}$
instead of $n$ and $f$) yields
\begin{equation}
\left(  f^{\circ a}\right)  ^{\circ\left(  m+1\right)  }=\left(  f^{\circ
a}\right)  ^{\circ\left(  \left(  m+1\right)  -1\right)  }\circ f^{\circ
a}=\left(  f^{\circ a}\right)  ^{\circ m}\circ f^{\circ a}
\label{pf.prop.ind.map-powers.2.7}%
\end{equation}
(since $\left(  m+1\right)  -1=m$).

But $a\left(  m+1\right)  =am+a$. Thus,%
\begin{align*}
f^{\circ\left(  a\left(  m+1\right)  \right)  }  &  =f^{\circ\left(
am+a\right)  }=\underbrace{f^{\circ\left(  am\right)  }}_{=\left(  f^{\circ
a}\right)  ^{\circ m}}\circ f^{\circ a}\\
&  \ \ \ \ \ \ \ \ \ \ \left(
\begin{array}
[c]{c}%
\text{by Proposition \ref{prop.ind.map-powers.2} \textbf{(a)}}\\
\text{(applied to }am\text{ and }a\text{ instead of }a\text{ and }b\text{)}%
\end{array}
\right) \\
&  =\left(  f^{\circ a}\right)  ^{\circ m}\circ f^{\circ a}=\left(  f^{\circ
a}\right)  ^{\circ\left(  m+1\right)  }\ \ \ \ \ \ \ \ \ \ \left(  \text{by
(\ref{pf.prop.ind.map-powers.2.7})}\right)  .
\end{align*}
In other words, (\ref{pf.prop.ind.map-powers.2.goal}) holds for $b=m+1$. This
completes the induction step. Thus, (\ref{pf.prop.ind.map-powers.2.goal}) is
proven by induction. Hence, Proposition \ref{prop.ind.map-powers.2}
\textbf{(b)} is proven.
\end{proof}

Note that Proposition \ref{prop.ind.map-powers.2} is similar to the rules of
exponents%
\[
n^{a+b}=n^{a}n^{b}\ \ \ \ \ \ \ \ \ \ \text{and}\ \ \ \ \ \ \ \ \ \ n^{ab}%
=\left(  n^{a}\right)  ^{b}%
\]
that hold for $n\in\mathbb{Q}$ and $a,b\in\mathbb{N}$ (and for various other
situations). Can we find similar analogues for other rules of exponents, such
as $\left(  mn\right)  ^{a}=m^{a}n^{a}$? The simplest analogue one could think
of for this rule would be $\left(  f\circ g\right)  ^{\circ a}=f^{\circ
a}\circ g^{\circ a}$; but this does not hold in general (unless $a\leq1$).
However, it turns out that this does hold if we assume that $f\circ g=g\circ
f$ (which is not automatically true, unlike the analogous equality $mn=nm$ for
integers). Let us prove this:

\begin{proposition}
\label{prop.ind.map-powers.3}Let $X$ be a set. Let $f:X\rightarrow X$ and
$g:X\rightarrow X$ be two maps such that $f\circ g=g\circ f$. Then:

\textbf{(a)} We have $f\circ g^{\circ b}=g^{\circ b}\circ f$ for each
$b\in\mathbb{N}$.

\textbf{(b)} We have $f^{\circ a}\circ g^{\circ b}=g^{\circ b}\circ f^{\circ
a}$ for each $a\in\mathbb{N}$ and $b\in\mathbb{N}$.

\textbf{(c)} We have $\left(  f\circ g\right)  ^{\circ a}=f^{\circ a}\circ
g^{\circ a}$ for each $a\in\mathbb{N}$.
\end{proposition}

\begin{example}
Let us see why the requirement $f\circ g=g\circ f$ is needed in Proposition
\ref{prop.ind.map-powers.3}:

Let $X$ be the set $\mathbb{Z}$. Let $f:X\rightarrow X$ be the map that sends
every integer $x$ to $-x$. Let $g:X\rightarrow X$ be the map that sends every
integer $x$ to $1-x$. Then, $f^{\circ2}=\operatorname*{id}\nolimits_{X}$
(since $f^{\circ2}\left(  x\right)  =f\left(  f\left(  x\right)  \right)
=-\left(  -x\right)  =x$ for each $x\in X$) and $g^{\circ2}=\operatorname*{id}%
\nolimits_{X}$ (since $g^{\circ2}\left(  x\right)  =g\left(  g\left(
x\right)  \right)  =1-\left(  1-x\right)  =x$ for each $x\in X$). But the map
$f\circ g$ satisfies $\left(  f\circ g\right)  \left(  x\right)  =f\left(
g\left(  x\right)  \right)  =-\left(  1-x\right)  =x-1$ for each $x\in X$.
Hence, $\left(  f\circ g\right)  ^{\circ2}\left(  x\right)  =\left(  f\circ
g\right)  \left(  \left(  f\circ g\right)  \left(  x\right)  \right)  =\left(
x-1\right)  -1=x-2$ for each $x\in X$. Thus, $\left(  f\circ g\right)
^{\circ2}\neq\operatorname*{id}\nolimits_{X}$. Comparing this with
$\underbrace{f^{\circ2}}_{=\operatorname*{id}\nolimits_{X}}\circ
\underbrace{g^{\circ2}}_{=\operatorname*{id}\nolimits_{X}}=\operatorname*{id}%
\nolimits_{X}\circ\operatorname*{id}\nolimits_{X}=\operatorname*{id}%
\nolimits_{X}$, we obtain $\left(  f\circ g\right)  ^{\circ2}\neq f^{\circ
2}\circ g^{\circ2}$. This shows that Proposition \ref{prop.ind.map-powers.3}
\textbf{(c)} would not hold without the requirement $f\circ g=g\circ f$.
\end{example}

\begin{proof}
[Proof of Proposition \ref{prop.ind.map-powers.3}.]\textbf{(a)} We claim that%
\begin{equation}
f\circ g^{\circ b}=g^{\circ b}\circ f\ \ \ \ \ \ \ \ \ \ \text{for each }%
b\in\mathbb{N}. \label{pf.prop.ind.map-powers.3.a.1}%
\end{equation}

Indeed, let us prove (\ref{pf.prop.ind.map-powers.3.a.1}) by induction on $b$:

\textit{Induction base:} We have $g^{\circ0}=\operatorname*{id}\nolimits_{X}$
(by (\ref{eq.def.ind.gen-ass-maps.fn.f0}), applied to $g$ instead of $f$).
Hence, $f\circ\underbrace{g^{\circ0}}_{=\operatorname*{id}\nolimits_{X}%
}=f\circ\operatorname*{id}\nolimits_{X}=f$ and $\underbrace{g^{\circ0}%
}_{=\operatorname*{id}\nolimits_{X}}\circ f=\operatorname*{id}\nolimits_{X}%
\circ f=f$. Comparing these two equalities, we obtain $f\circ g^{\circ
0}=g^{\circ0}\circ f$. In other words, (\ref{pf.prop.ind.map-powers.3.a.1})
holds for $b=0$. This completes the induction base.

\textit{Induction step:} Let $m\in\mathbb{N}$. Assume that
(\ref{pf.prop.ind.map-powers.3.a.1}) holds for $b=m$. We must prove that
(\ref{pf.prop.ind.map-powers.3.a.1}) holds for $b=m+1$.

We have assumed that (\ref{pf.prop.ind.map-powers.3.a.1}) holds for $b=m$. In
other words,%
\begin{equation}
f\circ g^{\circ m}=g^{\circ m}\circ f. \label{pf.prop.ind.map-powers.3.a.2}%
\end{equation}

Proposition \ref{prop.ind.gen-ass-maps.fgh} (applied to $Y=X$, $Z=X$, $W=X$,
$c=g$, $b=g^{\circ m}$ and $a=f$) yields
\begin{equation}
\left(  f\circ g^{\circ m}\right)  \circ g=f\circ\left(  g^{\circ m}\circ
g\right)  . \label{pf.prop.ind.map-powers.3.a.3}%
\end{equation}

Proposition \ref{prop.ind.map-powers.2} \textbf{(a)} (applied to $g$, $m$ and
$1$ instead of $f$, $a$ and $b$) yields
\begin{equation}
g^{\circ\left(  m+1\right)  }=g^{\circ m}\circ\underbrace{g^{\circ1}}%
_{=g}=g^{\circ m}\circ g. \label{pf.prop.ind.map-powers.3.a.4}%
\end{equation}
Hence,%
\begin{align}
f\circ\underbrace{g^{\circ\left(  m+1\right)  }}_{=g^{\circ m}\circ g}  &
=f\circ\left(  g^{\circ m}\circ g\right)  =\underbrace{\left(  f\circ g^{\circ
m}\right)  }_{\substack{=g^{\circ m}\circ f\\\text{(by
(\ref{pf.prop.ind.map-powers.3.a.2}))}}}\circ g\ \ \ \ \ \ \ \ \ \ \left(
\text{by (\ref{pf.prop.ind.map-powers.3.a.3})}\right) \nonumber\\
&  =\left(  g^{\circ m}\circ f\right)  \circ g=g^{\circ m}\circ
\underbrace{\left(  f\circ g\right)  }_{=g\circ f}\nonumber\\
&  \ \ \ \ \ \ \ \ \ \ \left(
\begin{array}
[c]{c}%
\text{by Proposition \ref{prop.ind.gen-ass-maps.fgh} (applied}\\
\text{to }Y=X\text{, }Z=X\text{, }W=X\text{, }c=g\text{, }b=f\text{ and
}a=g^{\circ m}\text{)}%
\end{array}
\right) \nonumber\\
&  =g^{\circ m}\circ\left(  g\circ f\right)  .
\label{pf.prop.ind.map-powers.3.a.5}%
\end{align}

On the other hand, Proposition \ref{prop.ind.gen-ass-maps.fgh} (applied to
$Y=X$, $Z=X$, $W=X$, $c=f$, $b=g$ and $a=g^{\circ m}$) yields
\[
\left(  g^{\circ m}\circ g\right)  \circ f=g^{\circ m}\circ\left(  g\circ
f\right)  .
\]
Comparing this with (\ref{pf.prop.ind.map-powers.3.a.5}), we obtain%
\[
f\circ g^{\circ\left(  m+1\right)  }=\underbrace{\left(  g^{\circ m}\circ
g\right)  }_{\substack{=g^{\circ\left(  m+1\right)  }\\\text{(by
(\ref{pf.prop.ind.map-powers.3.a.4}))}}}\circ f=g^{\circ\left(  m+1\right)
}\circ f.
\]
In other words, (\ref{pf.prop.ind.map-powers.3.a.1}) holds for $b=m+1$. This
completes the induction step. Thus, (\ref{pf.prop.ind.map-powers.3.a.1}) is
proven by induction.

Therefore, Proposition \ref{prop.ind.map-powers.3} \textbf{(a)} follows.

\textbf{(b)} Let $a\in\mathbb{N}$ and $b\in\mathbb{N}$. From $f\circ g=g\circ
f$, we obtain $g\circ f=f\circ g$. Hence, Proposition
\ref{prop.ind.map-powers.3} \textbf{(a)} (applied to $g$, $f$ and $a$ instead
of $f$, $g$ and $b$) yields $g\circ f^{\circ a}=f^{\circ a}\circ g$. In other
words, $f^{\circ a}\circ g=g\circ f^{\circ a}$. Hence, Proposition
\ref{prop.ind.map-powers.3} \textbf{(a)} (applied to $f^{\circ a}$ instead of
$f$) yields $f^{\circ a}\circ g^{\circ b}=g^{\circ b}\circ f^{\circ a}$. This
proves Proposition \ref{prop.ind.map-powers.3} \textbf{(b)}.

\textbf{(c)} We claim that%
\begin{equation}
\left(  f\circ g\right)  ^{\circ a}=f^{\circ a}\circ g^{\circ a}%
\ \ \ \ \ \ \ \ \ \ \text{for each }a\in\mathbb{N}.
\label{pf.prop.ind.map-powers.3.c.claim}%
\end{equation}

Indeed, let us prove (\ref{pf.prop.ind.map-powers.3.c.claim}) by induction on
$a$:

\textit{Induction base:} From (\ref{eq.def.ind.gen-ass-maps.fn.f0}), we obtain
$f^{\circ0}=\operatorname*{id}\nolimits_{X}$ and $g^{\circ0}%
=\operatorname*{id}\nolimits_{X}$ and $\left(  f\circ g\right)  ^{\circ
0}=\operatorname*{id}\nolimits_{X}$. Thus,%
\[
\left(  f\circ g\right)  ^{\circ0}=\operatorname*{id}\nolimits_{X}%
=\underbrace{\operatorname*{id}\nolimits_{X}}_{=f^{\circ0}}\circ
\underbrace{\operatorname*{id}\nolimits_{X}}_{=g^{\circ0}}=f^{\circ0}\circ
g^{\circ0}.
\]
In other words, (\ref{pf.prop.ind.map-powers.3.c.claim}) holds for $a=0$. This
completes the induction base.

\textit{Induction step:} Let $m\in\mathbb{N}$. Assume that
(\ref{pf.prop.ind.map-powers.3.c.claim}) holds for $a=m$. We must prove that
(\ref{pf.prop.ind.map-powers.3.c.claim}) holds for $a=m+1$.

We have assumed that (\ref{pf.prop.ind.map-powers.3.c.claim}) holds for $a=m$.
In other words,
\begin{equation}
\left(  f\circ g\right)  ^{\circ m}=f^{\circ m}\circ g^{\circ m}.
\label{pf.prop.ind.map-powers.3.c.0}%
\end{equation}

But Proposition \ref{prop.ind.map-powers.2} \textbf{(a)} (applied to $g$, $m$
and $1$ instead of $f$, $a$ and $b$) yields
\[
g^{\circ\left(  m+1\right)  }=g^{\circ m}\circ\underbrace{g^{\circ1}}%
_{=g}=g^{\circ m}\circ g.
\]
The same argument (applied to $f$ instead of $g$) yields $f^{\circ\left(
m+1\right)  }=f^{\circ m}\circ f$. Hence,%
\begin{equation}
\underbrace{f^{\circ\left(  m+1\right)  }}_{=f^{\circ m}\circ f}\circ
g^{\circ\left(  m+1\right)  }=\left(  f^{\circ m}\circ f\right)  \circ
g^{\circ\left(  m+1\right)  }=f^{\circ m}\circ\left(  f\circ g^{\circ\left(
m+1\right)  }\right)  \label{pf.prop.ind.map-powers.3.c.1}%
\end{equation}
(by Proposition \ref{prop.ind.gen-ass-maps.fgh} (applied to $Y=X$, $Z=X$,
$W=X$, $c=g^{\circ\left(  m+1\right)  }$, $b=f$ and $a=f^{\circ m}$)).

But Proposition \ref{prop.ind.map-powers.3} \textbf{(a)} (applied to $b=m+1$)
yields%
\[
f\circ g^{\circ\left(  m+1\right)  }=\underbrace{g^{\circ\left(  m+1\right)
}}_{=g^{\circ m}\circ g}\circ f=\left(  g^{\circ m}\circ g\right)  \circ
f=g^{\circ m}\circ\left(  g\circ f\right)
\]
(by Proposition \ref{prop.ind.gen-ass-maps.fgh} (applied to $Y=X$, $Z=X$,
$W=X$, $c=f$, $b=g$ and $a=g^{\circ m}$)). Hence,%
\begin{equation}
f\circ g^{\circ\left(  m+1\right)  }=g^{\circ m}\circ\underbrace{\left(
g\circ f\right)  }_{\substack{=f\circ g\\\text{(since }f\circ g=g\circ
f\text{)}}}=g^{\circ m}\circ\left(  f\circ g\right)  .
\label{pf.prop.ind.map-powers.3.c.3}%
\end{equation}

On the other hand, Proposition \ref{prop.ind.map-powers.2} \textbf{(a)}
(applied to $f\circ g$, $m$ and $1$ instead of $f$, $a$ and $b$) yields
\[
\left(  f\circ g\right)  ^{\circ\left(  m+1\right)  }=\underbrace{\left(
f\circ g\right)  ^{\circ m}}_{\substack{=f^{\circ m}\circ g^{\circ
m}\\\text{(by (\ref{pf.prop.ind.map-powers.3.c.0}))}}}\circ\underbrace{\left(
f\circ g\right)  ^{\circ1}}_{=f\circ g}=\left(  f^{\circ m}\circ g^{\circ
m}\right)  \circ\left(  f\circ g\right)  =f^{\circ m}\circ\left(  g^{\circ
m}\circ\left(  f\circ g\right)  \right)
\]
(by Proposition \ref{prop.ind.gen-ass-maps.fgh} (applied to $Y=X$, $Z=X$,
$W=X$, $c=f\circ g$, $b=g^{\circ m}$ and $a=f^{\circ m}$)). Hence,%
\[
\left(  f\circ g\right)  ^{\circ\left(  m+1\right)  }=f^{\circ m}%
\circ\underbrace{\left(  g^{\circ m}\circ\left(  f\circ g\right)  \right)
}_{\substack{=f\circ g^{\circ\left(  m+1\right)  }\\\text{(by
(\ref{pf.prop.ind.map-powers.3.c.3}))}}}=f^{\circ m}\circ\left(  f\circ
g^{\circ\left(  m+1\right)  }\right)  =f^{\circ\left(  m+1\right)  }\circ
g^{\circ\left(  m+1\right)  }%
\]
(by (\ref{pf.prop.ind.map-powers.3.c.1})). In other words,
(\ref{pf.prop.ind.map-powers.3.c.claim}) holds for $a=m+1$. This completes the
induction step. Thus, (\ref{pf.prop.ind.map-powers.3.c.claim}) is proven by
induction. Therefore, Proposition \ref{prop.ind.map-powers.3} \textbf{(c)} follows.
\end{proof}

\begin{remark}
In our above proof of Proposition \ref{prop.ind.map-powers.3}, we have not
used the notation $f_{n}\circ f_{n-1}\circ\cdots\circ f_{1}$ introduced in
Definition \ref{def.ind.gen-ass-maps.comp}, but instead relied on parentheses
and compositions of two maps (i.e., we have never composed more than two maps
at the same time). Thus, for example, in the proof of Proposition
\ref{prop.ind.map-powers.3} \textbf{(a)}, we wrote \textquotedblleft$\left(
g^{\circ m}\circ g\right)  \circ f$\textquotedblright\ and \textquotedblleft%
$g^{\circ m}\circ\left(  g\circ f\right)  $\textquotedblright\ rather than
\textquotedblleft$g^{\circ m}\circ g\circ f$\textquotedblright. But Remark
\ref{rmk.ind.gen-ass-maps.drop} says that we could have just as well dropped
all the parentheses. This would have saved us the trouble of explicitly
applying Proposition \ref{prop.ind.gen-ass-maps.fgh} (since if we drop all
parentheses, then there is no difference between \textquotedblleft$\left(
g^{\circ m}\circ g\right)  \circ f$\textquotedblright\ and \textquotedblleft%
$g^{\circ m}\circ\left(  g\circ f\right)  $\textquotedblright\ any more). This
way, the induction step in the proof of Proposition
\ref{prop.ind.map-powers.3} \textbf{(a)} could have been made much shorter:

\textit{Induction step (second version):} Let $m\in\mathbb{N}$. Assume that
(\ref{pf.prop.ind.map-powers.3.a.1}) holds for $b=m$. We must prove that
(\ref{pf.prop.ind.map-powers.3.a.1}) holds for $b=m+1$.

We have assumed that (\ref{pf.prop.ind.map-powers.3.a.1}) holds for $b=m$. In
other words,%
\begin{equation}
f\circ g^{\circ m}=g^{\circ m}\circ f.
\label{pf.prop.ind.map-powers.3.a.2nd.2}%
\end{equation}

Proposition \ref{prop.ind.map-powers.2} \textbf{(a)} (applied to $g$, $m$ and
$1$ instead of $f$, $a$ and $b$) yields $g^{\circ\left(  m+1\right)
}=g^{\circ m}\circ\underbrace{g^{\circ1}}_{=g}=g^{\circ m}\circ g$. Hence,%
\[
f\circ\underbrace{g^{\circ\left(  m+1\right)  }}_{=g^{\circ m}\circ
g}=\underbrace{f\circ g^{\circ m}}_{\substack{=g^{\circ m}\circ f\\\text{(by
(\ref{pf.prop.ind.map-powers.3.a.2nd.2}))}}}\circ g=g^{\circ m}\circ
\underbrace{f\circ g}_{=g\circ f}=\underbrace{g^{\circ m}\circ g}%
_{=g^{\circ\left(  m+1\right)  }}\circ f=g^{\circ\left(  m+1\right)  }\circ
f.
\]
In other words, (\ref{pf.prop.ind.map-powers.3.a.1}) holds for $b=m+1$. This
completes the induction step.

Similarly, we can simplify the proof of Proposition
\ref{prop.ind.map-powers.3} \textbf{(c)} by dropping the parentheses. (The
details are left to the reader.)
\end{remark}

\subsubsection{Composition of invertible maps}

The composition of two invertible maps is always invertible, and its inverse
can be computed by the following formula:

\begin{proposition}
\label{prop.ind.inverse-fg}Let $X$, $Y$ and $Z$ be three sets. Let
$b:X\rightarrow Y$ and $a:Y\rightarrow Z$ be two invertible maps. Then, the
map $a\circ b:X\rightarrow Z$ is invertible as well, and its inverse is%
\[
\left(  a\circ b\right)  ^{-1}=b^{-1}\circ a^{-1}.
\]

\end{proposition}

\begin{vershort}
\begin{proof}
[Proof of Proposition \ref{prop.ind.inverse-fg}.]As we have explained in
Definition \ref{def.ind.gen-ass-maps.comp}, we can drop the parentheses when
composing several maps. This will allow us to write expressions like
$b^{-1}\circ a^{-1}\circ a\circ b$ without specifying where parentheses should
be placed, and then pretending that they are placed wherever we would find
them most convenient.

The equalities%
\[
\left(  b^{-1}\circ a^{-1}\right)  \circ\left(  a\circ b\right)  =b^{-1}%
\circ\underbrace{a^{-1}\circ a}_{=\operatorname*{id}\nolimits_{Y}}\circ
b=b^{-1}\circ b=\operatorname*{id}\nolimits_{X}%
\]
and%
\[
\left(  a\circ b\right)  \circ\left(  b^{-1}\circ a^{-1}\right)
=a\circ\underbrace{b\circ b^{-1}}_{=\operatorname*{id}\nolimits_{Y}}\circ
a^{-1}=a\circ a^{-1}=\operatorname*{id}\nolimits_{Z}%
\]
show that the map $b^{-1}\circ a^{-1}$ is an inverse of $a\circ b$. Thus, the
map $a\circ b$ has an inverse (namely, $b^{-1}\circ a^{-1}$). In other words,
the map $a\circ b$ is invertible. Its inverse is $\left(  a\circ b\right)
^{-1}=b^{-1}\circ a^{-1}$ (since $b^{-1}\circ a^{-1}$ is an inverse of $a\circ
b$). This completes the proof of Proposition \ref{prop.ind.inverse-fg}.
\end{proof}
\end{vershort}

\begin{verlong}
\begin{proof}
[Proof of Proposition \ref{prop.ind.inverse-fg}.]Proposition
\ref{prop.ind.gen-ass-maps.fgh} (applied to $Y$, $b$, $a$ and $a^{-1}$ instead
of $W$, $c$, $b$ and $a$) shows that%
\[
\left(  a^{-1}\circ a\right)  \circ b=a^{-1}\circ\left(  a\circ b\right)  .
\]
Thus,%
\[
a^{-1}\circ\left(  a\circ b\right)  =\underbrace{\left(  a^{-1}\circ a\right)
}_{=\operatorname*{id}\nolimits_{Y}}\circ b=\operatorname*{id}\nolimits_{Y}%
\circ b=b.
\]

Proposition \ref{prop.ind.gen-ass-maps.fgh} (applied to $Z$, $Y$, $X$, $a\circ
b$, $a^{-1}$ and $b^{-1}$ instead of $Y$, $Z$, $W$, $c$, $b$ and $a$) shows
that%
\[
\left(  b^{-1}\circ a^{-1}\right)  \circ\left(  a\circ b\right)  =b^{-1}%
\circ\underbrace{\left(  a^{-1}\circ\left(  a\circ b\right)  \right)  }%
_{=b}=b^{-1}\circ b=\operatorname*{id}\nolimits_{X}.
\]

Proposition \ref{prop.ind.gen-ass-maps.fgh} (applied to $Z$, $Y$, $X$, $Y$,
$a^{-1}$, $b^{-1}$ and $b$ instead of $X$, $Y$, $Z$, $W$, $c$, $b$ and $a$)
shows that%
\[
\left(  b\circ b^{-1}\right)  \circ a^{-1}=b\circ\left(  b^{-1}\circ
a^{-1}\right)  .
\]
Thus,%
\[
b\circ\left(  b^{-1}\circ a^{-1}\right)  =\underbrace{\left(  b\circ
b^{-1}\right)  }_{=\operatorname*{id}\nolimits_{Y}}\circ a^{-1}%
=\operatorname*{id}\nolimits_{Y}\circ a^{-1}=a^{-1}.
\]

Proposition \ref{prop.ind.gen-ass-maps.fgh} (applied to $Z$, $X$, $Y$, $Z$,
$b^{-1}\circ a^{-1}$, $b$ and $a$ instead of $X$, $Y$, $Z$, $W$, $c$, $b$ and
$a$) shows that%
\[
\left(  a\circ b\right)  \circ\left(  b^{-1}\circ a^{-1}\right)
=a\circ\underbrace{\left(  b\circ\left(  b^{-1}\circ a^{-1}\right)  \right)
}_{=a^{-1}}=a\circ a^{-1}=\operatorname*{id}\nolimits_{Z}.
\]
Combining the equalities $\left(  b^{-1}\circ a^{-1}\right)  \circ\left(
a\circ b\right)  =\operatorname*{id}\nolimits_{X}$ and $\left(  a\circ
b\right)  \circ\left(  b^{-1}\circ a^{-1}\right)  =\operatorname*{id}%
\nolimits_{Z}$, we conclude that the map $b^{-1}\circ a^{-1}$ is an inverse of
$a\circ b$. Thus, the map $a\circ b$ has an inverse (namely, $b^{-1}\circ
a^{-1}$). In other words, the map $a\circ b$ is invertible. Its inverse is
$\left(  a\circ b\right)  ^{-1}=b^{-1}\circ a^{-1}$ (since $b^{-1}\circ
a^{-1}$ is an inverse of $a\circ b$). This completes the proof of Proposition
\ref{prop.ind.inverse-fg}.
\end{proof}
\end{verlong}

By a straightforward induction, we can extend Proposition
\ref{prop.ind.inverse-fg} to compositions of $n$ invertible maps:

\begin{proposition}
\label{prop.ind.inverse-comp}Let $n\in\mathbb{N}$. Let $X_{1},X_{2}%
,\ldots,X_{n+1}$ be $n+1$ sets. For each $i\in\left\{  1,2,\ldots,n\right\}
$, let $f_{i}:X_{i}\rightarrow X_{i+1}$ be an invertible map. Then, the map
$f_{n}\circ f_{n-1}\circ\cdots\circ f_{1}:X_{1}\rightarrow X_{n+1}$ is
invertible as well, and its inverse is%
\[
\left(  f_{n}\circ f_{n-1}\circ\cdots\circ f_{1}\right)  ^{-1}=f_{1}^{-1}\circ
f_{2}^{-1}\circ\cdots\circ f_{n}^{-1}.
\]

\end{proposition}

\begin{exercise}
\label{exe.ind.inverse-comp}Prove Proposition \ref{prop.ind.inverse-comp}.
\end{exercise}

In particular, Proposition \ref{prop.ind.inverse-comp} shows that any
composition of invertible maps is invertible. Since invertible maps are the
same as bijective maps, we can rewrite this as follows: Any composition of
bijective maps is bijective.

\subsection{\label{sect.ind.gen-com}General commutativity for addition of
numbers}

\subsubsection{The setup and the problem}

Throughout Section \ref{sect.ind.gen-com}, we let $\mathbb{A}$ be one of the
sets $\mathbb{N}$, $\mathbb{Z}$, $\mathbb{Q}$, $\mathbb{R}$ and $\mathbb{C}$.
The elements of $\mathbb{A}$ will be simply called \textit{numbers}.

There is an analogue of Proposition \ref{prop.ind.gen-ass-maps.fgh} for numbers:

\begin{proposition}
\label{prop.ind.gen-com.fgh}Let $a$, $b$ and $c$ be three numbers (i.e.,
elements of $\mathbb{A}$). Then, $\left(  a+b\right)  +c=a+\left(  b+c\right)
$.
\end{proposition}

Proposition \ref{prop.ind.gen-com.fgh} is known as the \textit{associativity
of addition} (in $\mathbb{A}$), and is fundamental; its proof can be found in
any textbook on the construction of the number system\footnote{For example,
Proposition \ref{prop.ind.gen-com.fgh} is proven in \cite[Theorem 3.2.3
(3)]{Swanso18} for the case when $\mathbb{A}=\mathbb{N}$; in \cite[Theorem
3.5.4 (3)]{Swanso18} for the case when $\mathbb{A}=\mathbb{Z}$; in
\cite[Theorem 3.6.4 (3)]{Swanso18} for the case when $\mathbb{A}=\mathbb{Q}$;
in \cite[Theorem 3.7.11]{Swanso18} for the case when $\mathbb{A}=\mathbb{R}$;
in \cite[Theorem 3.9.2]{Swanso18} for the case when $\mathbb{A}=\mathbb{C}$.}.

In Section \ref{sect.ind.gen-ass}, we have used Proposition
\ref{prop.ind.gen-ass-maps.fgh} to show that we can \textquotedblleft drop the
parentheses\textquotedblright\ in a composition $f_{n}\circ f_{n-1}\circ
\cdots\circ f_{1}$ of maps (i.e., all possible complete parenthesizations of
this composition are actually the same map). Likewise, we can use Proposition
\ref{prop.ind.gen-com.fgh} to show that we can \textquotedblleft drop the
parentheses\textquotedblright\ in a sum $a_{1}+a_{2}+\cdots+a_{n}$ of numbers
(i.e., all possible complete parenthesizations of this sum are actually the
same number). For example, if $a,b,c,d$ are four numbers, then the complete
parenthesizations of $a+b+c+d$ are%
\begin{align*}
&  \left(  \left(  a+b\right)  +c\right)  +d,\ \ \ \ \ \ \ \ \ \ \left(
a+\left(  b+c\right)  \right)  +d,\ \ \ \ \ \ \ \ \ \ \left(  a+b\right)
+\left(  c+d\right)  ,\\
&  a+\left(  \left(  b+c\right)  +d\right)  ,\ \ \ \ \ \ \ \ \ \ a+\left(
b+\left(  c+d\right)  \right)  ,
\end{align*}
and all of these five complete parenthesizations are the same number.

However, numbers behave better than maps. In particular, along with
Proposition \ref{prop.ind.gen-com.fgh}, they satisfy another law that maps
(generally) don't satisfy:

\begin{proposition}
\label{prop.ind.gen-com.fg}Let $a$ and $b$ be two numbers (i.e., elements of
$\mathbb{A}$). Then, $a+b=b+a$.
\end{proposition}

Proposition \ref{prop.ind.gen-com.fg} is known as the \textit{commutativity of
addition} (in $\mathbb{A}$), and again is a fundamental result whose proofs
are found in standard textbooks\footnote{For example, Proposition
\ref{prop.ind.gen-com.fg} is proven in \cite[Theorem 3.2.3 (4)]{Swanso18} for
the case when $\mathbb{A}=\mathbb{N}$; in \cite[Theorem 3.5.4 (4)]{Swanso18}
for the case when $\mathbb{A}=\mathbb{Z}$; in \cite[Theorem 3.6.4
(4)]{Swanso18} for the case when $\mathbb{A}=\mathbb{Q}$; in \cite[Theorem
3.7.11]{Swanso18} for the case when $\mathbb{A}=\mathbb{R}$; in \cite[Theorem
3.9.2]{Swanso18} for the case when $\mathbb{A}=\mathbb{C}$.}.

Furthermore, numbers can \textbf{always} be added, whereas maps can only be
composed if the domain of one is the codomain of the other. Thus, when we want
to take the sum of $n$ numbers $a_{1},a_{2},\ldots,a_{n}$, we can not only
choose where to put the parentheses, but also in what order the numbers should
appear in the sum. It turns out that neither of these choices affects the
result. For example, if $a,b,c$ are three numbers, then all $12$ possible sums%
\begin{align*}
&  \left(  a+b\right)  +c,\ \ \ \ \ \ \ \ \ \ a+\left(  b+c\right)
,\ \ \ \ \ \ \ \ \ \ \left(  a+c\right)  +b,\ \ \ \ \ \ \ \ \ \ a+\left(
c+b\right)  ,\\
&  \left(  b+a\right)  +c,\ \ \ \ \ \ \ \ \ \ b+\left(  a+c\right)
,\ \ \ \ \ \ \ \ \ \ \left(  b+c\right)  +a,\ \ \ \ \ \ \ \ \ \ b+\left(
c+a\right)  ,\\
&  \left(  c+a\right)  +b,\ \ \ \ \ \ \ \ \ \ c+\left(  a+b\right)
,\ \ \ \ \ \ \ \ \ \ \left(  c+b\right)  +a,\ \ \ \ \ \ \ \ \ \ c+\left(
b+a\right)
\end{align*}
are actually the same number. The reader can easily verify this for three
numbers $a,b,c$ (using Proposition \ref{prop.ind.gen-com.fgh} and Proposition
\ref{prop.ind.gen-com.fg}), but of course the general case (with $n$ numbers)
is more difficult. The independence of the result on the parenthesization can
be proven using the same arguments that we gave in Section
\ref{sect.ind.gen-ass} (except that the $\circ$ symbol is now replaced by
$+$), but the independence on the order cannot easily be shown (or even
stated) in this way.

Thus, we shall proceed differently: We shall rigorously define the sum of $n$
numbers without specifying an order in which they are added or using
parentheses. Unlike the composition of $n$ maps, which was defined for an
\textit{ordered list} of $n$ maps, we shall define the sum of $n$ numbers for
a \textit{family} of $n$ numbers (see the next subsection for the definition
of a \textquotedblleft family\textquotedblright). Families don't come with an
ordering chosen in advance, so we cannot single out any specific ordering for
use in the definition. Thus, the independence on the order will be baked right
into the definition.

Different solutions to the problem of formalizing the concept of the sum of
$n$ numbers can be found in \cite[Chapter 1, \S 1.5]{Bourba74}%
\footnote{Bourbaki, in \cite[Chapter 1, \S 1.5]{Bourba74}, define something
more general than a sum of $n$ numbers: They define the \textquotedblleft
composition\textquotedblright\ of a finite family of elements of a commutative
magma. The sum of $n$ numbers is a particular case of this concept when the
magma is the set $\mathbb{A}$ (endowed with its addition).}, in \cite[\S 3.3]%
{GalQua18} and in \cite[\S 2.4]{Clemen22}.

\subsubsection{Families}

Let us first define what we mean by a \textquotedblleft
family\textquotedblright\ of $n$ numbers. More generally, we can define a
family of elements of any set, or even a family of elements of
\textbf{different} sets. To motivate the definition, we first recall a concept
of an $n$-tuple:

\begin{remark}
\label{rmk.ind.families.tups}Let $n\in\mathbb{N}$.

\textbf{(a)} Let $A$ be a set. Then, to specify an $n$\textit{-tuple of
elements of }$A$ means specifying an element $a_{i}$ of $A$ for each
$i\in\left\{  1,2,\ldots,n\right\}  $. This $n$-tuple is then denoted by
$\left(  a_{1},a_{2},\ldots,a_{n}\right)  $ or by $\left(  a_{i}\right)
_{i\in\left\{  1,2,\ldots,n\right\}  }$. For each $i\in\left\{  1,2,\ldots
,n\right\}  $, we refer to $a_{i}$ as the $i$\textit{-th entry} of this $n$-tuple.

The set of all $n$-tuples of elements of $A$ is denoted by $A^{n}$ or by
$A^{\times n}$; it is called the $n$\textit{-th Cartesian power} of the set
$A$.

\textbf{(b)} More generally, we can define $n$-tuples of elements from
\textbf{different} sets: For each $i\in\left\{  1,2,\ldots,n\right\}  $, let
$A_{i}$ be a set. Then, to specify an $n$\textit{-tuple of elements of }%
$A_{1},A_{2},\ldots,A_{n}$ means specifying an element $a_{i}$ of $A_{i}$ for
each $i\in\left\{  1,2,\ldots,n\right\}  $. This $n$-tuple is (again) denoted
by $\left(  a_{1},a_{2},\ldots,a_{n}\right)  $ or by $\left(  a_{i}\right)
_{i\in\left\{  1,2,\ldots,n\right\}  }$. For each $i\in\left\{  1,2,\ldots
,n\right\}  $, we refer to $a_{i}$ as the $i$\textit{-th entry} of this $n$-tuple.

The set of all $n$-tuples of elements of $A_{1},A_{2},\ldots,A_{n}$ is denoted
by $A_{1}\times A_{2}\times\cdots\times A_{n}$ or by $\prod_{i=1}^{n}A_{i}$;
it is called the \textit{Cartesian product} of the $n$ sets $A_{1}%
,A_{2},\ldots,A_{n}$. These $n$ sets $A_{1},A_{2},\ldots,A_{n}$ are called the
\textit{factors} of this Cartesian product.
\end{remark}

\begin{example}
\label{exa.ind.families.tups}\textbf{(a)} The $3$-tuple $\left(  7,8,9\right)
$ is a $3$-tuple of elements of $\mathbb{N}$, and also a $3$-tuple of elements
of $\mathbb{Z}$. It can also be written in the form $\left(  6+i\right)
_{i\in\left\{  1,2,3\right\}  }$. Thus, $\left(  6+i\right)  _{i\in\left\{
1,2,3\right\}  }=\left(  6+1,6+2,6+3\right)  =\left(  7,8,9\right)
\in\mathbb{N}^{3}$ and also $\left(  6+i\right)  _{i\in\left\{  1,2,3\right\}
}\in\mathbb{Z}^{3}$.

\textbf{(b)} The $5$-tuple $\left(  \left\{  1\right\}  ,\left\{  2\right\}
,\left\{  3\right\}  ,\varnothing,\mathbb{N}\right)  $ is a $5$-tuple of
elements of the powerset of $\mathbb{N}$ (since $\left\{  1\right\}  ,\left\{
2\right\}  ,\left\{  3\right\}  ,\varnothing,\mathbb{N}$ are subsets of
$\mathbb{N}$, thus elements of the powerset of $\mathbb{N}$).

\textbf{(c)} The $0$-tuple $\left(  {}\right)  $ can be viewed as a $0$-tuple
of elements of \textbf{any} set $A$.

\textbf{(d)} If we let $\left[  n\right]  $ be the set $\left\{
1,2,\ldots,n\right\}  $ for each $n\in\mathbb{N}$, then $\left(
1,2,2,3,3\right)  $ is a $5$-tuple of elements of $\left[  1\right]  ,\left[
2\right]  ,\left[  3\right]  ,\left[  4\right]  ,\left[  5\right]  $ (because
$1\in\left[  1\right]  $, $2\in\left[  2\right]  $, $2\in\left[  3\right]  $,
$3\in\left[  4\right]  $ and $3\in\left[  5\right]  $). In other words,
$\left(  1,2,2,3,3\right)  \in\left[  1\right]  \times\left[  2\right]
\times\left[  3\right]  \times\left[  4\right]  \times\left[  5\right]  $.

\textbf{(e)} A $2$-tuple is the same as an ordered pair. A $3$-tuple is the
same as an ordered triple. A $1$-tuple of elements of a set $A$ is
\textquotedblleft almost\textquotedblright\ the same as a single element of
$A$; more precisely, there is a bijection
\[
A\rightarrow A^{1},\ \ \ \ \ \ \ \ \ \ a\mapsto\left(  a\right)
\]
from $A$ to the set of $1$-tuples of elements of $A$.
\end{example}

The notation \textquotedblleft$\left(  a_{i}\right)  _{i\in\left\{
1,2,\ldots,n\right\}  }$\textquotedblright\ in Remark
\ref{rmk.ind.families.tups} should be pronounced as \textquotedblleft the
$n$-tuple whose $i$-th entry is $a_{i}$ for each $i\in\left\{  1,2,\ldots
,n\right\}  $\textquotedblright. The letter \textquotedblleft$i$%
\textquotedblright\ is used as a variable in this notation (similar to the
\textquotedblleft$i$\textquotedblright\ in the expression \textquotedblleft%
$\sum_{i=1}^{n}i$\textquotedblright\ or in the expression \textquotedblleft
the map $\mathbb{N}\rightarrow\mathbb{N},\ i\mapsto i+1$\textquotedblright\ or
in the expression \textquotedblleft for all $i\in\mathbb{N}$, we have
$i+1>i$\textquotedblright); it does not refer to any specific element of
$\left\{  1,2,\ldots,n\right\}  $. As usual, it does not matter which letter
we are using for this variable (as long as it does not already have a
different meaning); thus, for example, the $3$-tuples $\left(  6+i\right)
_{i\in\left\{  1,2,3\right\}  }$ and $\left(  6+j\right)  _{j\in\left\{
1,2,3\right\}  }$ and $\left(  6+x\right)  _{x\in\left\{  1,2,3\right\}  }$
are all identical (and equal $\left(  7,8,9\right)  $).

We also note that the \textquotedblleft$\prod$\textquotedblright\ sign in
Remark \ref{rmk.ind.families.tups} \textbf{(b)} has a different meaning than
the \textquotedblleft$\prod$\textquotedblright\ sign in Section
\ref{sect.sums-repetitorium}. The former stands for a Cartesian product of
sets, whereas the latter stands for a product of numbers. In particular, a
product $\prod_{i=1}^{n}a_{i}$ of numbers does not change when its factors are
swapped, whereas a Cartesian product $\prod_{i=1}^{n}A_{i}$ of sets does. (In
particular, if $A$ and $B$ are two sets, then $A\times B$ and $B\times A$ are
different sets in general. The $2$-tuple $\left(  1,-1\right)  $ belongs to
$\mathbb{N}\times\mathbb{Z}$, but not to $\mathbb{Z}\times\mathbb{N}$.)

Thus, the purpose of an $n$-tuple is storing several elements (possibly of
different sets) in one \textquotedblleft container\textquotedblright. This is
a highly useful notion, but sometimes one wants a more general concept, which
can store several elements but not necessarily organized in a
\textquotedblleft linear order\textquotedblright. For example, assume you want
to store four integers $a,b,c,d$ in the form of a rectangular table $\left(
\begin{array}
[c]{cc}%
a & b\\
c & d
\end{array}
\right)  $ (also known as a \textquotedblleft$2\times2$-table of
integers\textquotedblright). Such a table doesn't have a well-defined
\textquotedblleft$1$-st entry\textquotedblright\ or \textquotedblleft$2$-nd
entry\textquotedblright\ (unless you agree on a specific order in which you
read it); instead, it makes sense to speak of a \textquotedblleft$\left(
1,2\right)  $-th entry\textquotedblright\ (i.e., the entry in row $1$ and
column $2$, which is $b$) or of a \textquotedblleft$\left(  2,2\right)  $-th
entry\textquotedblright\ (i.e., the entry in row $2$ and column $2$, which is
$d$). Thus, such tables work similarly to $n$-tuples, but they are
\textquotedblleft indexed\textquotedblright\ by pairs $\left(  i,j\right)  $
of appropriate integers rather than by the numbers $1,2,\ldots,n$.

The concept of a \textquotedblleft family\textquotedblright\ generalizes both
$n$-tuples and rectangular tables: It allows the entries to be indexed by the
elements of an arbitrary (possibly infinite) set $I$ instead of the numbers
$1,2,\ldots,n$. Here is its definition (watch the similarities to Remark
\ref{rmk.ind.families.tups}):

\begin{definition}
\label{def.ind.families.fams}Let $I$ be a set.

\textbf{(a)} Let $A$ be a set. Then, to specify an $I$\textit{-family of
elements of }$A$ means specifying an element $a_{i}$ of $A$ for each $i\in I$.
This $I$-family is then denoted by $\left(  a_{i}\right)  _{i\in I}$. For each
$i\in I$, we refer to $a_{i}$ as the $i$\textit{-th entry} of this $I$-family.
(Unlike the case of $n$-tuples, there is no notation like $\left(  a_{1}%
,a_{2},\ldots,a_{n}\right)  $ for $I$-families, because there is no natural
way in which their entries should be listed.)

An $I$-family of elements of $A$ is also called an $A$\textit{-valued }%
$I$\textit{-family}.

The set of all $I$-families of elements of $A$ is denoted by $A^{I}$ or by
$A^{\times I}$. (Note that the notation $A^{I}$ is also used for the set of
all maps from $I$ to $A$. But this set is more or less the same as the set of
all $I$-families of elements of $A$; see Remark \ref{rmk.ind.families.maps}
below for the details.)

\textbf{(b)} More generally, we can define $I$-families of elements from
\textbf{different} sets: For each $i\in I$, let $A_{i}$ be a set. Then, to
specify an $I$\textit{-family of elements of }$\left(  A_{i}\right)  _{i\in
I}$ means specifying an element $a_{i}$ of $A_{i}$ for each $i\in I$. This
$I$-family is (again) denoted by $\left(  a_{i}\right)  _{i\in I}$. For each
$i\in I$, we refer to $a_{i}$ as the $i$\textit{-th entry} of this $I$-family.

The set of all $I$-families of elements of $\left(  A_{i}\right)  _{i\in I}$
is denoted by $\prod_{i\in I}A_{i}$.

The word \textquotedblleft$I$-family\textquotedblright\ (without further
qualifications) means an $I$-family of elements of $\left(  A_{i}\right)
_{i\in I}$ for some sets $A_{i}$.

The word \textquotedblleft family\textquotedblright\ (without further
qualifications) means an $I$-family for some set $I$.
\end{definition}

\begin{example}
\label{exa.ind.families.fams}\textbf{(a)} The family $\left(  6+i\right)
_{i\in\left\{  0,3,5\right\}  }$ is a $\left\{  0,3,5\right\}  $-family of
elements of $\mathbb{N}$ (that is, an $\mathbb{N}$-valued $\left\{
0,3,5\right\}  $-family). It has three entries: Its $0$-th entry is $6+0=6$;
its $3$-rd entry is $6+3=9$; its $5$-th entry is $6+5=11$. Of course, this
family is also a $\left\{  0,3,5\right\}  $-family of elements of $\mathbb{Z}%
$. If we squint hard enough, we can pretend that this family is simply the
$3$-tuple $\left(  6,9,11\right)  $; but this is not advisable, and also does
not extend to situations in which there is no natural order on the set $I$.

\textbf{(b)} Let $X$ be the set $\left\{  \text{\textquotedblleft
cat\textquotedblright},\text{\textquotedblleft chicken\textquotedblright%
},\text{\textquotedblleft dog\textquotedblright}\right\}  $ consisting of
three words. Then, we can define an $X$-family $\left(  a_{i}\right)  _{i\in
X}$ of elements of $\mathbb{N}$ by setting%
\[
a_{\text{\textquotedblleft cat\textquotedblright}}%
=4,\ \ \ \ \ \ \ \ \ \ a_{\text{\textquotedblleft chicken\textquotedblright}%
}=2,\ \ \ \ \ \ \ \ \ \ a_{\text{\textquotedblleft dog\textquotedblright}}=4.
\]
This family has $3$ entries, which are $4$, $2$ and $4$; but there is no
natural order on the set $X$, so we cannot identify it with a $3$-tuple.

We can also rewrite this family as
\[
\left(  \text{the number of legs of a typical specimen of animal }i\right)
_{i\in X}.
\]
Of course, not every family will have a description like this; sometimes a
family is just a choice of elements without any underlying pattern.

\textbf{(c)} If $I$ is the empty set $\varnothing$, and if $A$ is any set,
then there is exactly one $I$-family of elements of $A$; namely, the
\textit{empty family}. Indeed, specifying such a family means specifying no
elements at all, and there is just one way to do that. We can denote the empty
family by $\left(  {}\right)  $, just like the empty $0$-tuple.

\textbf{(d)} The family $\left(  \left\vert i\right\vert \right)
_{i\in\mathbb{Z}}$ is a $\mathbb{Z}$-family of elements of $\mathbb{N}$
(because $\left\vert i\right\vert $ is an element of $\mathbb{N}$ for each
$i\in\mathbb{Z}$). It can also be regarded as a $\mathbb{Z}$-family of
elements of $\mathbb{Z}$.

\textbf{(e)} If $I$ is the set $\left\{  1,2,\ldots,n\right\}  $ for some
$n\in\mathbb{N}$, and if $A$ is any set, then an $I$-family $\left(
a_{i}\right)  _{i\in\left\{  1,2,\ldots,n\right\}  }$ of elements of $A$ is
the same as an $n$-tuple of elements of $A$. The same holds for families and
$n$-tuples of elements from different sets. Thus, any $n$ sets $A_{1}%
,A_{2},\ldots,A_{n}$ satisfy $\prod_{i\in\left\{  1,2,\ldots,n\right\}  }%
A_{i}=\prod_{i=1}^{n}A_{i}$.
\end{example}

The notation \textquotedblleft$\left(  a_{i}\right)  _{i\in I}$%
\textquotedblright\ in Definition \ref{def.ind.families.fams} should be
pronounced as \textquotedblleft the $I$-family whose $i$-th entry is $a_{i}$
for each $i\in I$\textquotedblright. The letter \textquotedblleft%
$i$\textquotedblright\ is used as a variable in this notation (similar to the
\textquotedblleft$i$\textquotedblright\ in the expression \textquotedblleft%
$\sum_{i=1}^{n}i$\textquotedblright); it does not refer to any specific
element of $I$. As usual, it does not matter which letter we are using for
this variable (as long as it does not already have a different meaning); thus,
for example, the $\mathbb{Z}$-families $\left(  \left\vert i\right\vert
\right)  _{i\in\mathbb{Z}}$ and $\left(  \left\vert p\right\vert \right)
_{p\in\mathbb{Z}}$ and $\left(  \left\vert w\right\vert \right)
_{w\in\mathbb{Z}}$ are all identical.

\begin{remark}
\label{rmk.ind.families.maps}Let $I$ and $A$ be two sets. What is the
difference between an $A$-valued $I$-family and a map from $I$ to $A$ ? Both
of these objects consist of a choice of an element of $A$ for each $i\in I$.

The main difference is terminological: e.g., when we speak of a family, the
elements of $A$ that constitute it are called its \textquotedblleft
entries\textquotedblright, whereas for a map they are called its
\textquotedblleft images\textquotedblright\ or \textquotedblleft
values\textquotedblright. Also, the notations for them are different: The
$A$-valued $I$-family $\left(  a_{i}\right)  _{i\in I}$ corresponds to the map
$I\rightarrow A,\ i\mapsto a_{i}$.

There is also another, subtler difference: A map from $I$ to $A$
\textquotedblleft knows\textquotedblright\ what the set $A$ is (so that, for
example, the maps $\mathbb{N}\rightarrow\mathbb{N},\ i\mapsto i$ and
$\mathbb{N}\rightarrow\mathbb{Z},\ i\mapsto i$ are considered different, even
though they map every element of $\mathbb{N}$ to the same value); but an
$A$-valued $I$-family does not \textquotedblleft know\textquotedblright\ what
the set $A$ is (so that, for example, the $\mathbb{N}$-valued $\mathbb{N}%
$-family $\left(  i\right)  _{i\in\mathbb{N}}$ is considered identical with
the $\mathbb{Z}$-valued $\mathbb{N}$-family $\left(  i\right)  _{i\in
\mathbb{N}}$). This matters occasionally when one wants to consider maps or
families for different sets simultaneously; it is not relevant if we just work
with $A$-valued $I$-families (or maps from $I$ to $A$) for two fixed sets $I$
and $A$. And either way, these conventions are not universal across the
mathematical literature; for some authors, maps from $I$ to $A$ do not
\textquotedblleft know\textquotedblright\ what $A$ is, whereas other authors
want families to \textquotedblleft know\textquotedblright\ this too.

What is certainly true, independently of any conventions, is the following
fact: If $I$ and $A$ are two sets, then the map%
\begin{align*}
\left\{  \text{maps from }I\text{ to }A\right\}   &  \rightarrow\left\{
A\text{-valued }I\text{-families}\right\}  ,\\
f  &  \mapsto\left(  f\left(  i\right)  \right)  _{i\in I}%
\end{align*}
is bijective. (Its inverse map sends every $A$-valued $I$-family $\left(
a_{i}\right)  _{i\in I}$ to the map $I\rightarrow A,\ i\mapsto a_{i}$.) Thus,
there is little harm in equating $\left\{  \text{maps from }I\text{ to
}A\right\}  $ with $\left\{  A\text{-valued }I\text{-families}\right\}  $.
\end{remark}

We already know from Example \ref{exa.ind.families.fams} \textbf{(e)} that
$n$-tuples are a particular case of families; the same holds for rectangular tables:

\begin{definition}
\label{def.ind.families.rectab}Let $A$ be a set. Let $n\in\mathbb{N}$ and
$m\in\mathbb{N}$. Then, an $n\times m$\textit{-table} of elements of $A$ means
an $A$-valued $\left\{  1,2,\ldots,n\right\}  \times\left\{  1,2,\ldots
,m\right\}  $-family. According to Remark \ref{rmk.ind.families.maps}, this is
tantamount to saying that an $n\times m$-table of elements of $A$ means a map
from $\left\{  1,2,\ldots,n\right\}  \times\left\{  1,2,\ldots,m\right\}  $ to
$A$, except for notational differences (such as referring to the elements that
constitute the $n\times m$-table as \textquotedblleft
entries\textquotedblright\ rather than \textquotedblleft
values\textquotedblright) and for the fact that an $n\times m$-table does not
\textquotedblleft know\textquotedblright\ $A$ (whereas a map would do).

In future chapters, we shall consider \textquotedblleft$n\times m$%
-matrices\textquotedblright, which are defined as maps from $\left\{
1,2,\ldots,n\right\}  \times\left\{  1,2,\ldots,m\right\}  $ to $A$ rather
than as $A$-valued $\left\{  1,2,\ldots,n\right\}  \times\left\{
1,2,\ldots,m\right\}  $-families. We shall keep using the same notations for
them as for $n\times m$-tables, but unlike $n\times m$-tables, they will
\textquotedblleft know\textquotedblright\ $A$ (that is, two $n\times
m$-matrices with the same entries but different sets $A$ will be considered
different). Anyway, this difference is minor.
\end{definition}

\begin{noncompile}
(The following has been removed, since I don't use these notations.)

\textbf{(b)} Let $C=\left(  c_{i,j}\right)  _{\left(  i,j\right)  \in\left\{
1,2,\ldots,n\right\}  \times\left\{  1,2,\ldots,m\right\}  }$ be an $n\times
m$-table of elements of $A$. Then, $C$ is often written as%
\[
\left(
\begin{array}
[c]{cccc}%
c_{1,1} & c_{1,2} & \cdots & c_{1,m}\\
c_{2,1} & c_{2,2} & \cdots & c_{2,m}\\
\vdots & \vdots & \ddots & \vdots\\
c_{n,1} & c_{n,2} & \cdots & c_{n,m}%
\end{array}
\right)
\]
(that is, as a rectangular table with $n$ rows and $m$ columns such that the
entry in the $i$-th row and the $j$-th column is $c_{i,j}$). For any
$p\in\left\{  1,2,\ldots,n\right\}  $, we denote the $1\times m$-table%
\[
\left(  c_{p,j}\right)  _{\left(  i,j\right)  \in\left\{  1\right\}
\times\left\{  1,2,\ldots,m\right\}  }=\left(
\begin{array}
[c]{cccc}%
c_{p,1} & c_{p,2} & \cdots & c_{p,m}%
\end{array}
\right)
\]
as the $p$\textit{-th row} of $C$. For any $q\in\left\{  1,2,\ldots,m\right\}
$, we denote the $n\times1$-table%
\[
\left(  c_{i,q}\right)  _{\left(  i,j\right)  \in\left\{  1,2,\ldots
,n\right\}  \times\left\{  1\right\}  }=\left(
\begin{array}
[c]{c}%
c_{1,q}\\
c_{2,q}\\
\vdots\\
c_{n,q}%
\end{array}
\right)
\]
as the $q$\textit{-th column} of $C$.
\end{noncompile}

\subsubsection{A desirable definition}

We now know what an $\mathbb{A}$-valued $S$-family is (for some set $S$): It
is just a way of choosing some element of $\mathbb{A}$ for each $s\in S$. When
this element is called $a_{s}$, the $S$-family is called $\left(
a_{s}\right)  _{s\in S}$.

We now want to define the sum of an $\mathbb{A}$-valued $S$-family $\left(
a_{s}\right)  _{s\in S}$ when the set $S$ is finite. Actually, we have already
seen a definition of this sum (which is called $\sum_{s\in S}a_{s}$) in
Section \ref{sect.sums-repetitorium}. The only problem with that definition is
that we don't know yet that it is legitimate. Let us nevertheless recall it
(rewriting it using the notion of an $\mathbb{A}$-valued $S$-family):

\begin{definition}
\label{def.ind.gen-com.defsum1}If $S$ is a finite set, and if $\left(
a_{s}\right)  _{s\in S}$ is an $\mathbb{A}$-valued $S$-family, then we want to
define the number $\sum_{s\in S}a_{s}$. We define this number by recursion on
$\left\vert S\right\vert $ as follows:

\begin{itemize}
\item If $\left\vert S\right\vert =0$, then $\sum_{s\in S}a_{s}$ is defined to
be $0$.

\item Let $n\in\mathbb{N}$. Assume that we have defined $\sum_{s\in S}a_{s}$
for every finite set $S$ with $\left\vert S\right\vert =n$ and any
$\mathbb{A}$-valued $S$-family $\left(  a_{s}\right)  _{s\in S}$. Now, if $S$
is a finite set with $\left\vert S\right\vert =n+1$, and if $\left(
a_{s}\right)  _{s\in S}$ is any $\mathbb{A}$-valued $S$-family, then
$\sum_{s\in S}a_{s}$ is defined by picking any $t\in S$ and setting%
\begin{equation}
\sum_{s\in S}a_{s}=a_{t}+\sum_{s\in S\setminus\left\{  t\right\}  }a_{s}.
\label{eq.def.ind.gen-com.defsum1.step}%
\end{equation}

\end{itemize}
\end{definition}

As we already observed in Section \ref{sect.sums-repetitorium}, it is not
obvious that this definition is legitimate: The right hand side of
(\ref{eq.def.ind.gen-com.defsum1.step}) is defined using a choice of $t$, but
we want our value of $\sum_{s\in S}a_{s}$ to depend only on $S$ and $\left(
a_{s}\right)  _{s\in S}$ (not on some arbitrarily chosen $t\in S$). Thus, we
cannot use this definition yet. Our main goal in this section is to prove that
it is indeed legitimate.

\subsubsection{The set of all possible sums}

There are two ways to approach this goal. One is to prove the legitimacy of
Definition \ref{def.ind.gen-com.defsum1} by strong induction on $\left\vert
S\right\vert $; the statement $\mathcal{A}\left(  n\right)  $ that we would be
proving for each $n\in\mathbb{N}$ here would be saying that Definition
\ref{def.ind.gen-com.defsum1} is legitimate for all finite sets $S$ satisfying
$\left\vert S\right\vert =n$. This is not hard, but conceptually confusing, as
it would require us to use Definition \ref{def.ind.gen-com.defsum1} for
\textbf{some} sets $S$ while its legitimacy for other sets $S$ is yet unproven.

We prefer to proceed in a different way: We shall first define a set
$\operatorname*{Sums}\left(  \left(  a_{s}\right)  _{s\in S}\right)  $ for any
$\mathbb{A}$-valued $S$-family $\left(  a_{s}\right)  _{s\in S}$; this set
shall consist (roughly speaking) of \textquotedblleft all possible values that
$\sum_{s\in S}a_{s}$ could have according to Definition
\ref{def.ind.gen-com.defsum1}\textquotedblright. This set will be defined
recursively, more or less following Definition \ref{def.ind.gen-com.defsum1},
but instead of relying on a choice of \textbf{some} $t\in S$, it will use
\textbf{all} possible elements $t\in S$. (See Definition
\ref{def.ind.gen-com.Sums} for the precise definition.) Unlike $\sum_{s\in
S}a_{s}$ itself, it will be a set of numbers, not a single number; however, it
has the advantage that the legitimacy of its definition will be immediately
obvious. Then, we will prove (in Theorem \ref{thm.ind.gen-com.Sums1}) that
this set $\operatorname*{Sums}\left(  \left(  a_{s}\right)  _{s\in S}\right)
$ is actually a $1$-element set; this will allow us to define $\sum_{s\in
S}a_{s}$ to be the unique element of $\operatorname*{Sums}\left(  \left(
a_{s}\right)  _{s\in S}\right)  $ for any $\mathbb{A}$-valued $S$-family
$\left(  a_{s}\right)  _{s\in S}$ (see Definition
\ref{def.ind.gen-com.defsum2}). Then, we will retroactively legitimize
Definition \ref{def.ind.gen-com.defsum1} by showing that Definition
\ref{def.ind.gen-com.defsum1} leads to the same value of $\sum_{s\in S}a_{s}$
as Definition \ref{def.ind.gen-com.defsum2} (no matter which $t\in S$ is
chosen). Having thus justified Definition \ref{def.ind.gen-com.defsum1}, we
will forget about the set $\operatorname*{Sums}\left(  \left(  a_{s}\right)
_{s\in S}\right)  $ and about Definition \ref{def.ind.gen-com.defsum2}.

In later subsections, we shall prove some basic properties of sums.

Let us define the set $\operatorname*{Sums}\left(  \left(  a_{s}\right)
_{s\in S}\right)  $, as promised:

\begin{definition}
\label{def.ind.gen-com.Sums}If $S$ is a finite set, and if $\left(
a_{s}\right)  _{s\in S}$ is an $\mathbb{A}$-valued $S$-family, then we want to
define the set $\operatorname*{Sums}\left(  \left(  a_{s}\right)  _{s\in
S}\right)  $ of numbers. We define this set by recursion on $\left\vert
S\right\vert $ as follows:

\begin{itemize}
\item If $\left\vert S\right\vert =0$, then $\operatorname*{Sums}\left(
\left(  a_{s}\right)  _{s\in S}\right)  $ is defined to be $\left\{
0\right\}  $.

\item Let $n\in\mathbb{N}$. Assume that we have defined $\operatorname*{Sums}%
\left(  \left(  a_{s}\right)  _{s\in S}\right)  $ for every finite set $S$
with $\left\vert S\right\vert =n$ and any $\mathbb{A}$-valued $S$-family
$\left(  a_{s}\right)  _{s\in S}$. Now, if $S$ is a finite set with
$\left\vert S\right\vert =n+1$, and if $\left(  a_{s}\right)  _{s\in S}$ is
any $\mathbb{A}$-valued $S$-family, then $\operatorname*{Sums}\left(  \left(
a_{s}\right)  _{s\in S}\right)  $ is defined by%
\begin{align}
&  \operatorname*{Sums}\left(  \left(  a_{s}\right)  _{s\in S}\right)
\nonumber\\
&  =\left\{  a_{t}+b\ \mid\ t\in S\text{ and }b\in\operatorname*{Sums}\left(
\left(  a_{s}\right)  _{s\in S\setminus\left\{  t\right\}  }\right)  \right\}
. \label{eq.def.ind.gen-com.Sums.rec}%
\end{align}
(The sets $\operatorname*{Sums}\left(  \left(  a_{s}\right)  _{s\in
S\setminus\left\{  t\right\}  }\right)  $ on the right hand side of this
equation are well-defined, because for each $t\in S$, we have $\left\vert
S\setminus\left\{  t\right\}  \right\vert =\left\vert S\right\vert -1=n$
(since $\left\vert S\right\vert =n+1$), and therefore $\operatorname*{Sums}%
\left(  \left(  a_{s}\right)  _{s\in S\setminus\left\{  t\right\}  }\right)  $
is well-defined by our assumption.)
\end{itemize}
\end{definition}

\begin{example}
\label{exa.def.ind.gen-com.Sums.1}Let $S$ be a finite set. Let $\left(
a_{s}\right)  _{s\in S}$ be an $\mathbb{A}$-valued $S$-family. Let us see what
Definition \ref{def.ind.gen-com.Sums} says when $S$ has only few elements:

\textbf{(a)} If $S=\varnothing$, then
\begin{equation}
\operatorname*{Sums}\left(  \left(  a_{s}\right)  _{s\in\varnothing}\right)
=\left\{  0\right\}  \label{eq.exa.def.ind.gen-com.Sums.1.0}%
\end{equation}
(directly by Definition \ref{def.ind.gen-com.Sums}, since $\left\vert
S\right\vert =\left\vert \varnothing\right\vert =0$ in this case).

\textbf{(b)} If $S=\left\{  x\right\}  $ for some element $x$, then Definition
\ref{def.ind.gen-com.Sums} yields%
\begin{align}
&  \operatorname*{Sums}\left(  \left(  a_{s}\right)  _{s\in\left\{  x\right\}
}\right) \nonumber\\
&  =\left\{  a_{t}+b\ \mid\ t\in\left\{  x\right\}  \text{ and }%
b\in\operatorname*{Sums}\left(  \left(  a_{s}\right)  _{s\in\left\{
x\right\}  \setminus\left\{  t\right\}  }\right)  \right\} \nonumber\\
&  =\left\{  a_{x}+b\ \mid\ b\in\operatorname*{Sums}\left(  \left(
a_{s}\right)  _{s\in\left\{  x\right\}  \setminus\left\{  x\right\}  }\right)
\right\}  \ \ \ \ \ \ \ \ \ \ \left(  \text{since the only }t\in\left\{
x\right\}  \text{ is }x\right) \nonumber\\
&  =\left\{  a_{x}+b\ \mid\ b\in\underbrace{\operatorname*{Sums}\left(
\left(  a_{s}\right)  _{s\in\varnothing}\right)  }_{=\left\{  0\right\}
}\right\}  \ \ \ \ \ \ \ \ \ \ \left(  \text{since }\left\{  x\right\}
\setminus\left\{  x\right\}  =\varnothing\right) \nonumber\\
&  =\left\{  a_{x}+b\ \mid\ b\in\left\{  0\right\}  \right\}  =\left\{
a_{x}+0\right\}  =\left\{  a_{x}\right\}  .
\label{eq.exa.def.ind.gen-com.Sums.1.1}%
\end{align}

\textbf{(c)} If $S=\left\{  x,y\right\}  $ for two distinct elements $x$ and
$y$, then Definition \ref{def.ind.gen-com.Sums} yields%
\begin{align*}
&  \operatorname*{Sums}\left(  \left(  a_{s}\right)  _{s\in\left\{
x,y\right\}  }\right) \\
&  =\left\{  a_{t}+b\ \mid\ t\in\left\{  x,y\right\}  \text{ and }%
b\in\operatorname*{Sums}\left(  \left(  a_{s}\right)  _{s\in\left\{
x,y\right\}  \setminus\left\{  t\right\}  }\right)  \right\} \\
&  =\left\{  a_{x}+b\ \mid\ b\in\operatorname*{Sums}\left(  \left(
a_{s}\right)  _{s\in\left\{  x,y\right\}  \setminus\left\{  x\right\}
}\right)  \right\} \\
&  \ \ \ \ \ \ \ \ \ \ \cup\left\{  a_{y}+b\ \mid\ b\in\operatorname*{Sums}%
\left(  \left(  a_{s}\right)  _{s\in\left\{  x,y\right\}  \setminus\left\{
y\right\}  }\right)  \right\} \\
&  =\left\{  a_{x}+b\ \mid\ b\in\underbrace{\operatorname*{Sums}\left(
\left(  a_{s}\right)  _{s\in\left\{  y\right\}  }\right)  }%
_{\substack{=\left\{  a_{y}\right\}  \\\text{(by
(\ref{eq.exa.def.ind.gen-com.Sums.1.1}), applied to }y\text{ instead of
}x\text{)}}}\right\} \\
&  \ \ \ \ \ \ \ \ \ \ \cup\left\{  a_{y}+b\ \mid\ b\in
\underbrace{\operatorname*{Sums}\left(  \left(  a_{s}\right)  _{s\in\left\{
x\right\}  }\right)  }_{\substack{=\left\{  a_{x}\right\}  \\\text{(by
(\ref{eq.exa.def.ind.gen-com.Sums.1.1}))}}}\right\} \\
&  \ \ \ \ \ \ \ \ \ \ \left(  \text{since }\left\{  x,y\right\}
\setminus\left\{  x\right\}  =\left\{  y\right\}  \text{ and }\left\{
x,y\right\}  \setminus\left\{  y\right\}  =\left\{  x\right\}  \right) \\
&  =\underbrace{\left\{  a_{x}+b\ \mid\ b\in\left\{  a_{y}\right\}  \right\}
}_{=\left\{  a_{x}+a_{y}\right\}  }\cup\underbrace{\left\{  a_{y}%
+b\ \mid\ b\in\left\{  a_{x}\right\}  \right\}  }_{=\left\{  a_{y}%
+a_{x}\right\}  }\\
&  =\left\{  a_{x}+a_{y}\right\}  \cup\left\{  a_{y}+a_{x}\right\}  =\left\{
a_{x}+a_{y},a_{y}+a_{x}\right\}  =\left\{  a_{x}+a_{y}\right\}
\end{align*}
(since $a_{y}+a_{x}=a_{x}+a_{y}$).

\textbf{(d)} Similar reasoning shows that if $S=\left\{  x,y,z\right\}  $ for
three distinct elements $x$, $y$ and $z$, then%
\[
\operatorname*{Sums}\left(  \left(  a_{s}\right)  _{s\in\left\{
x,y,z\right\}  }\right)  =\left\{  a_{x}+\left(  a_{y}+a_{z}\right)
,a_{y}+\left(  a_{x}+a_{z}\right)  ,a_{z}+\left(  a_{x}+a_{y}\right)
\right\}  .
\]
It is not hard to check (using Proposition \ref{prop.ind.gen-com.fgh} and
Proposition \ref{prop.ind.gen-com.fg}) that the three elements $a_{x}+\left(
a_{y}+a_{z}\right)  $, $a_{y}+\left(  a_{x}+a_{z}\right)  $ and $a_{z}+\left(
a_{x}+a_{y}\right)  $ of this set are equal, so we may call them $a_{x}%
+a_{y}+a_{z}$; thus, we can rewrite this equality as%
\[
\operatorname*{Sums}\left(  \left(  a_{s}\right)  _{s\in\left\{
x,y,z\right\}  }\right)  =\left\{  a_{x}+a_{y}+a_{z}\right\}  .
\]

\textbf{(e)} Going further, we can see that if $S=\left\{  x,y,z,w\right\}  $
for four distinct elements $x$, $y$, $z$ and $w$, then%
\begin{align*}
\operatorname*{Sums}\left(  \left(  a_{s}\right)  _{s\in\left\{
x,y,z,w\right\}  }\right)   &  =\left\{  a_{x}+\left(  a_{y}+a_{z}%
+a_{w}\right)  ,a_{y}+\left(  a_{x}+a_{z}+a_{w}\right)  ,\right. \\
&  \ \ \ \ \ \ \ \ \ \ \left.  a_{z}+\left(  a_{x}+a_{y}+a_{w}\right)
,a_{w}+\left(  a_{x}+a_{y}+a_{z}\right)  \right\}  .
\end{align*}
Again, it is not hard to prove that
\begin{align*}
a_{x}+\left(  a_{y}+a_{z}+a_{w}\right)   &  =a_{y}+\left(  a_{x}+a_{z}%
+a_{w}\right) \\
&  =a_{z}+\left(  a_{x}+a_{y}+a_{w}\right)  =a_{w}+\left(  a_{x}+a_{y}%
+a_{z}\right)  ,
\end{align*}
and thus the set $\operatorname*{Sums}\left(  \left(  a_{s}\right)
_{s\in\left\{  x,y,z,w\right\}  }\right)  $ is again a $1$-element set, whose
unique element can be called $a_{x}+a_{y}+a_{z}+a_{w}$.
\end{example}

These examples suggest that the set $\operatorname*{Sums}\left(  \left(
a_{s}\right)  _{s\in S}\right)  $ should always be a $1$-element set. This is
precisely what we are going to claim now:

\begin{theorem}
\label{thm.ind.gen-com.Sums1}If $S$ is a finite set, and if $\left(
a_{s}\right)  _{s\in S}$ is an $\mathbb{A}$-valued $S$-family, then the set
$\operatorname*{Sums}\left(  \left(  a_{s}\right)  _{s\in S}\right)  $ is a
$1$-element set.
\end{theorem}

\subsubsection{The set of all possible sums is a $1$-element set: proof}

Before we step to the proof of Theorem \ref{thm.ind.gen-com.Sums1}, we observe
an almost trivial lemma:

\begin{lemma}
\label{lem.ind.gen-com.Sums-lem}Let $a$, $b$ and $c$ be three numbers (i.e.,
elements of $\mathbb{A}$). Then, $a+\left(  b+c\right)  =b+\left(  a+c\right)
$.
\end{lemma}

\begin{proof}
[Proof of Lemma \ref{lem.ind.gen-com.Sums-lem}.]Proposition
\ref{prop.ind.gen-com.fgh} (applied to $b$ and $a$ instead of $a$ and $b$)
yields $\left(  b+a\right)  +c=b+\left(  a+c\right)  $. Also, Proposition
\ref{prop.ind.gen-com.fgh} yields $\left(  a+b\right)  +c=a+\left(
b+c\right)  $. Hence,%
\[
a+\left(  b+c\right)  =\underbrace{\left(  a+b\right)  }%
_{\substack{=b+a\\\text{(by Proposition \ref{prop.ind.gen-com.fg})}%
}}+c=\left(  b+a\right)  +c=b+\left(  a+c\right)  .
\]
This proves Lemma \ref{lem.ind.gen-com.Sums-lem}.
\end{proof}

\begin{proof}
[Proof of Theorem \ref{thm.ind.gen-com.Sums1}.]We shall prove Theorem
\ref{thm.ind.gen-com.Sums1} by strong induction on $\left\vert S\right\vert $:

Let $m\in\mathbb{N}$. Assume that Theorem \ref{thm.ind.gen-com.Sums1} holds
under the condition that $\left\vert S\right\vert <m$. We must now prove that
Theorem \ref{thm.ind.gen-com.Sums1} holds under the condition that $\left\vert
S\right\vert =m$.

We have assumed that Theorem \ref{thm.ind.gen-com.Sums1} holds under the
condition that $\left\vert S\right\vert <m$. In other words, the following
claim holds:

\begin{statement}
\textit{Claim 1:} Let $S$ be a finite set satisfying $\left\vert S\right\vert
<m$. Let $\left(  a_{s}\right)  _{s\in S}$ be an $\mathbb{A}$-valued
$S$-family. Then, the set $\operatorname*{Sums}\left(  \left(  a_{s}\right)
_{s\in S}\right)  $ is a $1$-element set.
\end{statement}

Now, we must now prove that Theorem \ref{thm.ind.gen-com.Sums1} holds under
the condition that $\left\vert S\right\vert =m$. In other words, we must prove
the following claim:

\begin{statement}
\textit{Claim 2:} Let $S$ be a finite set satisfying $\left\vert S\right\vert
=m$. Let $\left(  a_{s}\right)  _{s\in S}$ be an $\mathbb{A}$-valued
$S$-family. Then, the set $\operatorname*{Sums}\left(  \left(  a_{s}\right)
_{s\in S}\right)  $ is a $1$-element set.
\end{statement}

Before we start proving Claim 2, let us prove two auxiliary claims:

\begin{statement}
\textit{Claim 3:} Let $S$ be a finite set satisfying $\left\vert S\right\vert
<m$. Let $\left(  a_{s}\right)  _{s\in S}$ be an $\mathbb{A}$-valued
$S$-family. Let $r\in S$. Let $g\in\operatorname*{Sums}\left(  \left(
a_{s}\right)  _{s\in S}\right)  $ and $c\in\operatorname*{Sums}\left(  \left(
a_{s}\right)  _{s\in S\setminus\left\{  r\right\}  }\right)  $. Then,
$g=a_{r}+c$.
\end{statement}

[\textit{Proof of Claim 3:} The set $S\setminus\left\{  r\right\}  $ is a
subset of the finite set $S$, and thus itself is finite. Moreover, $r\in S$,
so that $\left\vert S\setminus\left\{  r\right\}  \right\vert =\left\vert
S\right\vert -1$. Thus, $\left\vert S\right\vert =\left\vert S\setminus
\left\{  r\right\}  \right\vert +1$. Hence, the definition of
$\operatorname*{Sums}\left(  \left(  a_{s}\right)  _{s\in S}\right)  $ yields%
\begin{equation}
\operatorname*{Sums}\left(  \left(  a_{s}\right)  _{s\in S}\right)  =\left\{
a_{t}+b\ \mid\ t\in S\text{ and }b\in\operatorname*{Sums}\left(  \left(
a_{s}\right)  _{s\in S\setminus\left\{  t\right\}  }\right)  \right\}  .
\label{pf.thm.ind.gen-com.Sums1.c3.pf.1}%
\end{equation}
But recall that $r\in S$ and $c\in\operatorname*{Sums}\left(  \left(
a_{s}\right)  _{s\in S\setminus\left\{  r\right\}  }\right)  $. Hence, the
number $a_{r}+c$ has the form $a_{t}+b$ for some $t\in S$ and $b\in
\operatorname*{Sums}\left(  \left(  a_{s}\right)  _{s\in S\setminus\left\{
t\right\}  }\right)  $ (namely, for $t=r$ and $b=c$). In other words,%
\[
a_{r}+c\in\left\{  a_{t}+b\ \mid\ t\in S\text{ and }b\in\operatorname*{Sums}%
\left(  \left(  a_{s}\right)  _{s\in S\setminus\left\{  t\right\}  }\right)
\right\}  .
\]
In view of (\ref{pf.thm.ind.gen-com.Sums1.c3.pf.1}), this rewrites as
$a_{r}+c\in\operatorname*{Sums}\left(  \left(  a_{s}\right)  _{s\in S}\right)
$.

But Claim 1 shows that the set $\operatorname*{Sums}\left(  \left(
a_{s}\right)  _{s\in S}\right)  $ is a $1$-element set. Hence, any two
elements of $\operatorname*{Sums}\left(  \left(  a_{s}\right)  _{s\in
S}\right)  $ are equal. In other words, any $x\in\operatorname*{Sums}\left(
\left(  a_{s}\right)  _{s\in S}\right)  $ and $y\in\operatorname*{Sums}\left(
\left(  a_{s}\right)  _{s\in S}\right)  $ satisfy $x=y$. Applying this to
$x=g$ and $y=a_{r}+c$, we obtain $g=a_{r}+c$ (since $g\in\operatorname*{Sums}%
\left(  \left(  a_{s}\right)  _{s\in S}\right)  $ and $a_{r}+c\in
\operatorname*{Sums}\left(  \left(  a_{s}\right)  _{s\in S}\right)  $). This
proves Claim 3.]

\begin{statement}
\textit{Claim 4:} Let $S$ be a finite set satisfying $\left\vert S\right\vert
=m$. Let $\left(  a_{s}\right)  _{s\in S}$ be an $\mathbb{A}$-valued
$S$-family. Let $p\in S$ and $q\in S$. Let $f\in\operatorname*{Sums}\left(
\left(  a_{s}\right)  _{s\in S\setminus\left\{  p\right\}  }\right)  $ and
$g\in\operatorname*{Sums}\left(  \left(  a_{s}\right)  _{s\in S\setminus
\left\{  q\right\}  }\right)  $. Then, $a_{p}+f=a_{q}+g$.
\end{statement}

[\textit{Proof of Claim 4:} We have $p\in S$, and thus $\left\vert
S\setminus\left\{  p\right\}  \right\vert =\left\vert S\right\vert
-1<\left\vert S\right\vert =m$. Hence, Claim 1 (applied to $S\setminus\left\{
p\right\}  $ instead of $S$) yields that the set $\operatorname*{Sums}\left(
\left(  a_{s}\right)  _{s\in S\setminus\left\{  p\right\}  }\right)  $ is a
$1$-element set. In other words, $\operatorname*{Sums}\left(  \left(
a_{s}\right)  _{s\in S\setminus\left\{  p\right\}  }\right)  $ can be written
in the form $\operatorname*{Sums}\left(  \left(  a_{s}\right)  _{s\in
S\setminus\left\{  p\right\}  }\right)  =\left\{  h\right\}  $ for some number
$h$. Consider this $h$.

We are in one of the following two cases:

\textit{Case 1:} We have $p=q$.

\textit{Case 2:} We have $p\neq q$.

Let us first consider Case 1. In this case, we have $p=q$. Hence, $q=p$, so
that $a_{q}=a_{p}$.

We have $f\in\operatorname*{Sums}\left(  \left(  a_{s}\right)  _{s\in
S\setminus\left\{  p\right\}  }\right)  =\left\{  h\right\}  $, so that $f=h$.
Also,
\begin{align*}
g  &  \in\operatorname*{Sums}\left(  \left(  a_{s}\right)  _{s\in
S\setminus\left\{  q\right\}  }\right)  =\operatorname*{Sums}\left(  \left(
a_{s}\right)  _{s\in S\setminus\left\{  p\right\}  }\right)
\ \ \ \ \ \ \ \ \ \ \left(  \text{since }q=p\right) \\
&  =\left\{  h\right\}  ,
\end{align*}
so that $g=h$. Comparing $a_{p}+\underbrace{f}_{=h}=a_{p}+h$ with
$\underbrace{a_{q}}_{=a_{p}}+\underbrace{g}_{=h}=a_{p}+h$, we obtain
$a_{p}+f=a_{q}+g$. Hence, Claim 4 is proven in Case 1.

Let us now consider Case 2. In this case, we have $p\neq q$. Thus, $q\neq p$,
so that $q\notin\left\{  p\right\}  $.

We have $S\setminus\left\{  p,q\right\}  \subseteq S\setminus\left\{
p\right\}  $ (since $\left\{  p\right\}  \subseteq\left\{  p,q\right\}  $) and
thus $\left\vert S\setminus\left\{  p,q\right\}  \right\vert \leq\left\vert
S\setminus\left\{  p\right\}  \right\vert <m$. Hence, Claim 1 (applied to
$S\setminus\left\{  p,q\right\}  $ instead of $S$) shows that the set
\newline$\operatorname*{Sums}\left(  \left(  a_{s}\right)  _{s\in
S\setminus\left\{  p,q\right\}  }\right)  $ is a $1$-element set. In other
words, $\operatorname*{Sums}\left(  \left(  a_{s}\right)  _{s\in
S\setminus\left\{  p,q\right\}  }\right)  $ has the form%
\[
\operatorname*{Sums}\left(  \left(  a_{s}\right)  _{s\in S\setminus\left\{
p,q\right\}  }\right)  =\left\{  c\right\}
\]
for some number $c$. Consider this $c$. Hence,%
\[
c\in\left\{  c\right\}  =\operatorname*{Sums}\left(  \left(  a_{s}\right)
_{s\in S\setminus\left\{  p,q\right\}  }\right)  =\operatorname*{Sums}\left(
\left(  a_{s}\right)  _{s\in\left(  S\setminus\left\{  p\right\}  \right)
\setminus\left\{  q\right\}  }\right)
\]
(since $S\setminus\left\{  p,q\right\}  =\left(  S\setminus\left\{  p\right\}
\right)  \setminus\left\{  q\right\}  $). Also, $q\in S\setminus\left\{
p\right\}  $ (since $q\in S$ and $q\notin\left\{  p\right\}  $). Thus, Claim 3
(applied to $S\setminus\left\{  p\right\}  $, $q$ and $f$ instead of $S$, $r$
and $g$) yields $f=a_{q}+c$ (since $\left\vert S\setminus\left\{  p\right\}
\right\vert <m$ and $f\in\operatorname*{Sums}\left(  \left(  a_{s}\right)
_{s\in S\setminus\left\{  p\right\}  }\right)  $ and $c\in\operatorname*{Sums}%
\left(  \left(  a_{s}\right)  _{s\in\left(  S\setminus\left\{  p\right\}
\right)  \setminus\left\{  q\right\}  }\right)  $).

\begin{vershort}
The same argument (but with $p$, $q$, $f$ and $g$ replaced by $q$, $p$, $g$
and $f$) yields $g=a_{p}+c$.
\end{vershort}

\begin{verlong}
Also, $q\in S$, and thus $\left\vert S\setminus\left\{  q\right\}  \right\vert
=\left\vert S\right\vert -1<\left\vert S\right\vert =m$. Furthermore, $p\neq
q$, so that $p\notin\left\{  q\right\}  $. Now,%
\[
c\in\operatorname*{Sums}\left(  \left(  a_{s}\right)  _{s\in S\setminus
\left\{  p,q\right\}  }\right)  =\operatorname*{Sums}\left(  \left(
a_{s}\right)  _{s\in\left(  S\setminus\left\{  q\right\}  \right)
\setminus\left\{  p\right\}  }\right)
\]
(since $S\setminus\left\{  p,q\right\}  =\left(  S\setminus\left\{  q\right\}
\right)  \setminus\left\{  p\right\}  $). Also, $p\in S\setminus\left\{
q\right\}  $ (since $p\in S$ and $p\notin\left\{  q\right\}  $). Thus, Claim 3
(applied to $S\setminus\left\{  q\right\}  $ and $p$ instead of $S$ and $r$)
yields $g=a_{p}+c$ (since $\left\vert S\setminus\left\{  q\right\}
\right\vert <m$ and $g\in\operatorname*{Sums}\left(  \left(  a_{s}\right)
_{s\in S\setminus\left\{  q\right\}  }\right)  $ and $c\in\operatorname*{Sums}%
\left(  \left(  a_{s}\right)  _{s\in\left(  S\setminus\left\{  q\right\}
\right)  \setminus\left\{  p\right\}  }\right)  $).
\end{verlong}

Now,%
\[
a_{p}+\underbrace{f}_{=a_{q}+c}=a_{p}+\left(  a_{q}+c\right)  =a_{q}+\left(
a_{p}+c\right)
\]
(by Lemma \ref{lem.ind.gen-com.Sums-lem}, applied to $a=a_{p}$ and $b=a_{q}$).
Comparing this with%
\[
a_{q}+\underbrace{g}_{=a_{p}+c}=a_{q}+\left(  a_{p}+c\right)  ,
\]
we obtain $a_{p}+f=a_{q}+g$. Thus, Claim 4 is proven in Case 2.

We have now proven Claim 4 in both Cases 1 and 2. Since these two Cases cover
all possibilities, we thus conclude that Claim 4 always holds.]

We can now prove Claim 2:

[\textit{Proof of Claim 2:} If $\left\vert S\right\vert =0$, then Claim 2
holds\footnote{\textit{Proof.} Assume that $\left\vert S\right\vert =0$.
Hence, Definition \ref{def.ind.gen-com.Sums} yields $\operatorname*{Sums}%
\left(  \left(  a_{s}\right)  _{s\in S}\right)  =\left\{  0\right\}  $. Hence,
the set $\operatorname*{Sums}\left(  \left(  a_{s}\right)  _{s\in S}\right)  $
is a $1$-element set (since the set $\left\{  0\right\}  $ is a $1$-element
set). In other words, Claim 2 holds. Qed.}. Hence, for the rest of this proof
of Claim 2, we can WLOG assume that we don't have $\left\vert S\right\vert
=0$. Assume this.

We have $\left\vert S\right\vert \neq0$ (since we don't have $\left\vert
S\right\vert =0$). Hence, $\left\vert S\right\vert $ is a positive integer.
Thus, $\left\vert S\right\vert -1\in\mathbb{N}$. Also, the set $S$ is nonempty
(since $\left\vert S\right\vert \neq0$). Hence, there exists some $p\in S$.
Consider this $p$.

We have $p\in S$ and thus $\left\vert S\setminus\left\{  p\right\}
\right\vert =\left\vert S\right\vert -1<\left\vert S\right\vert =m$. Hence,
Claim 1 (applied to $S\setminus\left\{  p\right\}  $ instead of $S$) shows
that the set $\operatorname*{Sums}\left(  \left(  a_{s}\right)  _{s\in
S\setminus\left\{  p\right\}  }\right)  $ is a $1$-element set. In other
words, $\operatorname*{Sums}\left(  \left(  a_{s}\right)  _{s\in
S\setminus\left\{  p\right\}  }\right)  $ has the form%
\[
\operatorname*{Sums}\left(  \left(  a_{s}\right)  _{s\in S\setminus\left\{
p\right\}  }\right)  =\left\{  f\right\}
\]
for some number $f$. Consider this $f$. Thus,%
\begin{equation}
f\in\left\{  f\right\}  =\operatorname*{Sums}\left(  \left(  a_{s}\right)
_{s\in S\setminus\left\{  p\right\}  }\right)  .
\label{pf.thm.ind.gen-com.Sums1.2}%
\end{equation}

Define $n\in\mathbb{N}$ by $n=\left\vert S\right\vert -1$. (This is allowed,
since $\left\vert S\right\vert -1\in\mathbb{N}$.) Then, from $n=\left\vert
S\right\vert -1$, we obtain $\left\vert S\right\vert =n+1$. Hence, the
definition of $\operatorname*{Sums}\left(  \left(  a_{s}\right)  _{s\in
S}\right)  $ yields%
\begin{equation}
\operatorname*{Sums}\left(  \left(  a_{s}\right)  _{s\in S}\right)  =\left\{
a_{t}+b\ \mid\ t\in S\text{ and }b\in\operatorname*{Sums}\left(  \left(
a_{s}\right)  _{s\in S\setminus\left\{  t\right\}  }\right)  \right\}  .
\label{pf.thm.ind.gen-com.Sums1.3}%
\end{equation}
But recall that $p\in S$ and $f\in\operatorname*{Sums}\left(  \left(
a_{s}\right)  _{s\in S\setminus\left\{  p\right\}  }\right)  $. Hence, the
number $a_{p}+f$ has the form $a_{t}+b$ for some $t\in S$ and $b\in
\operatorname*{Sums}\left(  \left(  a_{s}\right)  _{s\in S\setminus\left\{
t\right\}  }\right)  $ (namely, for $t=p$ and $b=f$). In other words,%
\[
a_{p}+f\in\left\{  a_{t}+b\ \mid\ t\in S\text{ and }b\in\operatorname*{Sums}%
\left(  \left(  a_{s}\right)  _{s\in S\setminus\left\{  t\right\}  }\right)
\right\}  .
\]
In view of (\ref{pf.thm.ind.gen-com.Sums1.3}), this rewrites as
\[
a_{p}+f\in\operatorname*{Sums}\left(  \left(  a_{s}\right)  _{s\in S}\right)
.
\]
Thus,%
\begin{equation}
\left\{  a_{p}+f\right\}  \subseteq\operatorname*{Sums}\left(  \left(
a_{s}\right)  _{s\in S}\right)  . \label{pf.thm.ind.gen-com.Sums1.5}%
\end{equation}

Next, we shall show the reverse inclusion (i.e., we shall show that
$\operatorname*{Sums}\left(  \left(  a_{s}\right)  _{s\in S}\right)
\subseteq\left\{  a_{p}+f\right\}  $).

Indeed, let $w\in\operatorname*{Sums}\left(  \left(  a_{s}\right)  _{s\in
S}\right)  $. Thus,%
\[
w\in\operatorname*{Sums}\left(  \left(  a_{s}\right)  _{s\in S}\right)
=\left\{  a_{t}+b\ \mid\ t\in S\text{ and }b\in\operatorname*{Sums}\left(
\left(  a_{s}\right)  _{s\in S\setminus\left\{  t\right\}  }\right)  \right\}
\]
(by (\ref{pf.thm.ind.gen-com.Sums1.3})). In other words, $w$ can be written as
$w=a_{t}+b$ for some $t\in S$ and $b\in\operatorname*{Sums}\left(  \left(
a_{s}\right)  _{s\in S\setminus\left\{  t\right\}  }\right)  $. Consider these
$t$ and $b$, and denote them by $q$ and $g$. Thus, $q\in S$ and $g\in
\operatorname*{Sums}\left(  \left(  a_{s}\right)  _{s\in S\setminus\left\{
q\right\}  }\right)  $ satisfy $w=a_{q}+g$.

But Claim 4 yields $a_{p}+f=a_{q}+g$. Comparing this with $w=a_{q}+g$, we
obtain $w=a_{p}+f$. Thus, $w\in\left\{  a_{p}+f\right\}  $.

Now, forget that we fixed $w$. We thus have proven that $w\in\left\{
a_{p}+f\right\}  $ for each $w\in\operatorname*{Sums}\left(  \left(
a_{s}\right)  _{s\in S}\right)  $. In other words, $\operatorname*{Sums}%
\left(  \left(  a_{s}\right)  _{s\in S}\right)  \subseteq\left\{
a_{p}+f\right\}  $. Combining this with (\ref{pf.thm.ind.gen-com.Sums1.5}), we
conclude that $\operatorname*{Sums}\left(  \left(  a_{s}\right)  _{s\in
S}\right)  =\left\{  a_{p}+f\right\}  $. Hence, the set $\operatorname*{Sums}%
\left(  \left(  a_{s}\right)  _{s\in S}\right)  $ is a $1$-element set. This
proves Claim 2.]

Now, we have proven Claim 2. But Claim 2 says precisely that Theorem
\ref{thm.ind.gen-com.Sums1} holds under the condition that $\left\vert
S\right\vert =m$. Hence, we have proven that Theorem
\ref{thm.ind.gen-com.Sums1} holds under the condition that $\left\vert
S\right\vert =m$. This completes the induction step. Thus, Theorem
\ref{thm.ind.gen-com.Sums1} is proven by strong induction.
\end{proof}

\subsubsection{Sums of numbers are well-defined}

We can now give a new definition of the sum $\sum_{s\in S}a_{s}$ (for any
finite set $S$ and any $\mathbb{A}$-valued $S$-family $\left(  a_{s}\right)
_{s\in S}$), which is different from Definition \ref{def.ind.gen-com.defsum1}
and (unlike the latter) is clearly legitimate:

\begin{definition}
\label{def.ind.gen-com.defsum2}Let $S$ be a finite set, and let $\left(
a_{s}\right)  _{s\in S}$ be an $\mathbb{A}$-valued $S$-family. Then, the set
$\operatorname*{Sums}\left(  \left(  a_{s}\right)  _{s\in S}\right)  $ is a
$1$-element set (by Theorem \ref{thm.ind.gen-com.Sums1}). We define
$\sum_{s\in S}a_{s}$ to be the unique element of this set
$\operatorname*{Sums}\left(  \left(  a_{s}\right)  _{s\in S}\right)  $.
\end{definition}

However, we have not reached our goal yet: After all, we wanted to prove that
Definition \ref{def.ind.gen-com.defsum1} is legitimate, rather than replace it
by a new definition!

Fortunately, we are very close to achieving this goal (after having done all
the hard work in the proof of Theorem \ref{thm.ind.gen-com.Sums1} above); we
are soon going to show that Definition \ref{def.ind.gen-com.defsum1} is
justified \textbf{and} that it is equivalent to Definition
\ref{def.ind.gen-com.defsum2} (that is, both definitions yield the same value
of $\sum_{s\in S}a_{s}$). First, we need a simple lemma, which says that the
notation $\sum_{s\in S}a_{s}$ defined in Definition
\ref{def.ind.gen-com.defsum2} \textquotedblleft behaves\textquotedblright%
\ like the one defined in Definition \ref{def.ind.gen-com.defsum1}:

\begin{lemma}
\label{lem.ind.gen-com.same-sum}In this lemma, we shall use Definition
\ref{def.ind.gen-com.defsum2} (not Definition \ref{def.ind.gen-com.defsum1}).

Let $S$ be a finite set, and let $\left(  a_{s}\right)  _{s\in S}$ be an
$\mathbb{A}$-valued $S$-family.

\textbf{(a)} If $\left\vert S\right\vert =0$, then
\[
\sum_{s\in S}a_{s}=0.
\]

\textbf{(b)} For any $t\in S$, we have%
\[
\sum_{s\in S}a_{s}=a_{t}+\sum_{s\in S\setminus\left\{  t\right\}  }a_{s}.
\]

\end{lemma}

\begin{proof}
[Proof of Lemma \ref{lem.ind.gen-com.same-sum}.]\textbf{(a)} Assume that
$\left\vert S\right\vert =0$. Thus, $\operatorname*{Sums}\left(  \left(
a_{s}\right)  _{s\in S}\right)  =\left\{  0\right\}  $ (by the definition of
$\operatorname*{Sums}\left(  \left(  a_{s}\right)  _{s\in S}\right)  $).
Hence, the unique element of the set $\operatorname*{Sums}\left(  \left(
a_{s}\right)  _{s\in S}\right)  $ is $0$.

But Definition \ref{def.ind.gen-com.defsum2} yields that $\sum_{s\in S}a_{s}$
is the unique element of the set $\operatorname*{Sums}\left(  \left(
a_{s}\right)  _{s\in S}\right)  $. Thus, $\sum_{s\in S}a_{s}$ is $0$ (since
the unique element of the set $\operatorname*{Sums}\left(  \left(
a_{s}\right)  _{s\in S}\right)  $ is $0$). In other words, $\sum_{s\in S}%
a_{s}=0$. This proves Lemma \ref{lem.ind.gen-com.same-sum} \textbf{(a)}.

\textbf{(b)} Let $p\in S$. Thus, $\left\vert S\setminus\left\{  p\right\}
\right\vert =\left\vert S\right\vert -1$, so that $\left\vert S\right\vert
=\left\vert S\setminus\left\{  p\right\}  \right\vert +1$. Hence, the
definition of $\operatorname*{Sums}\left(  \left(  a_{s}\right)  _{s\in
S}\right)  $ yields%
\begin{equation}
\operatorname*{Sums}\left(  \left(  a_{s}\right)  _{s\in S}\right)  =\left\{
a_{t}+b\ \mid\ t\in S\text{ and }b\in\operatorname*{Sums}\left(  \left(
a_{s}\right)  _{s\in S\setminus\left\{  t\right\}  }\right)  \right\}  .
\label{pf.lem.ind.gen-com.same-sum.b.1}%
\end{equation}

Definition \ref{def.ind.gen-com.defsum2} yields that $\sum_{s\in
S\setminus\left\{  p\right\}  }a_{s}$ is the unique element of the set
\newline$\operatorname*{Sums}\left(  \left(  a_{s}\right)  _{s\in
S\setminus\left\{  p\right\}  }\right)  $. Thus, $\sum_{s\in S\setminus
\left\{  p\right\}  }a_{s}\in\operatorname*{Sums}\left(  \left(  a_{s}\right)
_{s\in S\setminus\left\{  p\right\}  }\right)  $. Thus, $a_{p}+\sum_{s\in
S\setminus\left\{  p\right\}  }a_{s}$ is a number of the form $a_{t}+b$ for
some $t\in S$ and some $b\in\operatorname*{Sums}\left(  \left(  a_{s}\right)
_{s\in S\setminus\left\{  t\right\}  }\right)  $ (namely, for $t=p$ and
$b=\sum_{s\in S\setminus\left\{  p\right\}  }a_{s}$). In other words,%
\[
a_{p}+\sum_{s\in S\setminus\left\{  p\right\}  }a_{s}\in\left\{  a_{t}%
+b\ \mid\ t\in S\text{ and }b\in\operatorname*{Sums}\left(  \left(
a_{s}\right)  _{s\in S\setminus\left\{  t\right\}  }\right)  \right\}  .
\]
In view of (\ref{pf.lem.ind.gen-com.same-sum.b.1}), this rewrites as%
\begin{equation}
a_{p}+\sum_{s\in S\setminus\left\{  p\right\}  }a_{s}\in\operatorname*{Sums}%
\left(  \left(  a_{s}\right)  _{s\in S}\right)  .
\label{pf.lem.ind.gen-com.same-sum.b.3}%
\end{equation}

But Definition \ref{def.ind.gen-com.defsum2} yields that $\sum_{s\in S}a_{s}$
is the unique element of the set $\operatorname*{Sums}\left(  \left(
a_{s}\right)  _{s\in S}\right)  $. Hence, the set $\operatorname*{Sums}\left(
\left(  a_{s}\right)  _{s\in S}\right)  $ consists only of the element
$\sum_{s\in S}a_{s}$. In other words,
\[
\operatorname*{Sums}\left(  \left(  a_{s}\right)  _{s\in S}\right)  =\left\{
\sum_{s\in S}a_{s}\right\}  .
\]
Thus, (\ref{pf.lem.ind.gen-com.same-sum.b.3}) rewrites as
\[
a_{p}+\sum_{s\in S\setminus\left\{  p\right\}  }a_{s}\in\left\{  \sum_{s\in
S}a_{s}\right\}  .
\]
In other words, $a_{p}+\sum_{s\in S\setminus\left\{  p\right\}  }a_{s}%
=\sum_{s\in S}a_{s}$. Thus, $\sum_{s\in S}a_{s}=a_{p}+\sum_{s\in
S\setminus\left\{  p\right\}  }a_{s}$.

Now, forget that we fixed $p$. We thus have proven that for any $p\in S$, we
have $\sum_{s\in S}a_{s}=a_{p}+\sum_{s\in S\setminus\left\{  p\right\}  }%
a_{s}$. Renaming the variable $p$ as $t$ in this statement, we obtain the
following: For any $t\in S$, we have $\sum_{s\in S}a_{s}=a_{t}+\sum_{s\in
S\setminus\left\{  t\right\}  }a_{s}$. This proves Lemma
\ref{lem.ind.gen-com.same-sum} \textbf{(b)}.
\end{proof}

We can now finally state what we wanted to state:

\begin{theorem}
\label{thm.ind.gen-com.wd}\textbf{(a)} Definition
\ref{def.ind.gen-com.defsum1} is legitimate: i.e., the value of $\sum_{s\in
S}a_{s}$ in Definition \ref{def.ind.gen-com.defsum1} does not depend on the
choice of $t$.

\textbf{(b)} Definition \ref{def.ind.gen-com.defsum1} is equivalent to
Definition \ref{def.ind.gen-com.defsum2}: i.e., both of these definitions
yield the same value of $\sum_{s\in S}a_{s}$.
\end{theorem}

It makes sense to call Theorem \ref{thm.ind.gen-com.wd} \textbf{(a)} the
\textit{general commutativity theorem}, as it says that a sum of $n$ numbers
can be computed in an arbitrary order.

\begin{proof}
[Proof of Theorem \ref{thm.ind.gen-com.wd}.]Let us first use Definition
\ref{def.ind.gen-com.defsum2} (not Definition \ref{def.ind.gen-com.defsum1}).
Then, for any finite set $S$ and any $\mathbb{A}$-valued $S$-family $\left(
a_{s}\right)  _{s\in S}$, we can compute the number $\sum_{s\in S}a_{s}$ by
the following algorithm (which uses recursion on $\left\vert S\right\vert $):

\begin{itemize}
\item If $\left\vert S\right\vert =0$, then $\sum_{s\in S}a_{s}=0$. (This
follows from Lemma \ref{lem.ind.gen-com.same-sum} \textbf{(a)}.)

\item Otherwise, we have $\left\vert S\right\vert =n+1$ for some
$n\in\mathbb{N}$. Consider this $n$. Thus, $\left\vert S\right\vert
=n+1\geq1>0$, so that the set $S$ is nonempty. Fix any $t\in S$. (Such a $t$
exists, since the set $S$ is nonempty.) We have $\left\vert S\setminus\left\{
t\right\}  \right\vert =\left\vert S\right\vert -1=n$ (since $\left\vert
S\right\vert =n+1$), so that we can assume (because we are using recursion)
that $\sum_{s\in S\setminus\left\{  t\right\}  }a_{s}$ has already been
computed. Then, $\sum_{s\in S}a_{s}=a_{t}+\sum_{s\in S\setminus\left\{
t\right\}  }a_{s}$. (This follows from Lemma \ref{lem.ind.gen-com.same-sum}
\textbf{(b)}.)
\end{itemize}

We can restate this algorithm as an alternative definition of $\sum_{s\in
S}a_{s}$; it then takes the following form:

\begin{statement}
\textit{Alternative definition of }$\sum_{s\in S}a_{s}$\textit{ for any finite
set }$S$ \textit{and any }$\mathbb{A}$\textit{-valued }$S$\textit{-family
}$\left(  a_{s}\right)  _{s\in S}$\textit{:} If $S$ is a finite set, and if
$\left(  a_{s}\right)  _{s\in S}$ is an $\mathbb{A}$-valued $S$-family, then
we define $\sum_{s\in S}a_{s}$ by recursion on $\left\vert S\right\vert $ as follows:

\begin{itemize}
\item If $\left\vert S\right\vert =0$, then $\sum_{s\in S}a_{s}$ is defined to
be $0$.

\item Let $n\in\mathbb{N}$. Assume that we have defined $\sum_{s\in S}a_{s}$
for every finite set $S$ with $\left\vert S\right\vert =n$ and any
$\mathbb{A}$-valued $S$-family $\left(  a_{s}\right)  _{s\in S}$. Now, if $S$
is a finite set with $\left\vert S\right\vert =n+1$, and if $\left(
a_{s}\right)  _{s\in S}$ is any $\mathbb{A}$-valued $S$-family, then
$\sum_{s\in S}a_{s}$ is defined by picking any $t\in S$ and setting%
\begin{equation}
\sum_{s\in S}a_{s}=a_{t}+\sum_{s\in S\setminus\left\{  t\right\}  }a_{s}.
\label{pf.thm.ind.gen-com.wd.ad.t}%
\end{equation}

\end{itemize}
\end{statement}

This alternative definition of $\sum_{s\in S}a_{s}$ merely follows the above
algorithm for computing $\sum_{s\in S}a_{s}$. Thus, it is guaranteed to always
yield the same value of $\sum_{s\in S}a_{s}$ as Definition
\ref{def.ind.gen-com.defsum2}, independently of the choice of $t$. Hence, we
obtain the following:

\begin{statement}
\textit{Claim 1:} This alternative definition is legitimate (i.e., the value
of $\sum_{s\in S}a_{s}$ in (\ref{pf.thm.ind.gen-com.wd.ad.t}) does not depend
on the choice of $t$), and is equivalent to Definition
\ref{def.ind.gen-com.defsum2}.
\end{statement}

But on the other hand, this alternative definition is precisely Definition
\ref{def.ind.gen-com.defsum1}. Hence, Claim 1 rewrites as follows: Definition
\ref{def.ind.gen-com.defsum1} is legitimate (i.e., the value of $\sum_{s\in
S}a_{s}$ in Definition \ref{def.ind.gen-com.defsum1} does not depend on the
choice of $t$), and is equivalent to Definition \ref{def.ind.gen-com.defsum2}.
This proves both parts \textbf{(a)} and \textbf{(b)} of Theorem
\ref{thm.ind.gen-com.wd}.
\end{proof}

Theorem \ref{thm.ind.gen-com.wd} \textbf{(a)} shows that Definition
\ref{def.ind.gen-com.defsum1} is legitimate.

Thus, at last, we have vindicated the notation $\sum_{s\in S}a_{s}$ that was
introduced in Section \ref{sect.sums-repetitorium} (because the definition of
this notation we gave in Section \ref{sect.sums-repetitorium} was precisely
Definition \ref{def.ind.gen-com.defsum1}). We can now forget about Definition
\ref{def.ind.gen-com.defsum2}, since it has served its purpose (which was to
justify Definition \ref{def.ind.gen-com.defsum1}). (Of course, we could also
forget about Definition \ref{def.ind.gen-com.defsum1} instead, and use
Definition \ref{def.ind.gen-com.defsum2} as our definition of $\sum_{s\in
S}a_{s}$ (after all, these two definitions are equivalent, as we now know).
Then, we would have to replace every reference to the definition of
$\sum_{s\in S}a_{s}$ by a reference to Lemma \ref{lem.ind.gen-com.same-sum};
in particular, we would have to replace every use of (\ref{eq.sum.def.1}) by a
use of Lemma \ref{lem.ind.gen-com.same-sum} \textbf{(b)}. Other than this,
everything would work the same way.)

The notation $\sum_{s\in S}a_{s}$ has several properties, many of which were
collected in Section \ref{sect.sums-repetitorium}. We shall prove some of
these properties later in this section.

From now on, we shall be using all the conventions and notations regarding
sums that we introduced in Section \ref{sect.sums-repetitorium}. In
particular, expressions of the form \textquotedblleft$\sum_{s\in S}a_{s}%
+b$\textquotedblright\ shall always be interpreted as $\left(  \sum_{s\in
S}a_{s}\right)  +b$, not as $\sum_{s\in S}\left(  a_{s}+b\right)  $; but
expressions of the form \textquotedblleft$\sum_{s\in S}ba_{s}c$%
\textquotedblright\ shall always be understood to mean $\sum_{s\in S}\left(
ba_{s}c\right)  $.

\subsubsection{Triangular numbers revisited}

Recall one specific notation we introduced in Section
\ref{sect.sums-repetitorium}: If $u$ and $v$ are two integers, and if $a_{s}$
is a number for each $s\in\left\{  u,u+1,\ldots,v\right\}  $, then $\sum
_{s=u}^{v}a_{s}$ is defined by%
\[
\sum_{s=u}^{v}a_{s}=\sum_{s\in\left\{  u,u+1,\ldots,v\right\}  }a_{s}.
\]
This sum $\sum_{s=u}^{v}a_{s}$ is also denoted by $a_{u}+a_{u+1}+\cdots+a_{v}$.

We are now ready to do something that we evaded in Section
\ref{sect.ind.trinum}: namely, to speak of the sum of the first $n$ positive
integers without having to define it recursively. Indeed, we can now interpret
this sum as $\sum_{i\in\left\{  1,2,\ldots,n\right\}  }i$, an expression which
has a well-defined meaning because we have shown that the notation $\sum_{s\in
S}a_{s}$ is well-defined. We can also rewrite this expression as $\sum
_{i=1}^{n}i$ or as $1+2+\cdots+n$.

Thus, the classical fact that the sum of the first $n$ positive integers is
$\dfrac{n\left(  n+1\right)  }{2}$ can now be stated as follows:

\begin{proposition}
\label{prop.ind.gen-com.n(n+1)/2}We have%
\begin{equation}
\sum_{i\in\left\{  1,2,\ldots,n\right\}  }i=\dfrac{n\left(  n+1\right)  }%
{2}\ \ \ \ \ \ \ \ \ \ \text{for each }n\in\mathbb{N}.
\label{eq.prop.ind.gen-com.n(n+1)/2.claim}%
\end{equation}

\end{proposition}

\begin{proof}
[Proof of Proposition \ref{prop.ind.gen-com.n(n+1)/2}.]We shall prove
(\ref{eq.prop.ind.gen-com.n(n+1)/2.claim}) by induction on $n$:

\textit{Induction base:} We have $\left\{  1,2,\ldots,0\right\}  =\varnothing$
and thus $\left\vert \left\{  1,2,\ldots,0\right\}  \right\vert =\left\vert
\varnothing\right\vert =0$. Hence, the definition of $\sum_{i\in\left\{
1,2,\ldots,0\right\}  }i$ yields%
\begin{equation}
\sum_{i\in\left\{  1,2,\ldots,0\right\}  }i=0.
\label{pf.prop.ind.gen-com.n(n+1)/2.IB.1}%
\end{equation}
(To be more precise, we have used the first bullet point of Definition
\ref{def.ind.gen-com.defsum1} here, which says that $\sum_{s\in S}a_{s}=0$
whenever the set $S$ and the $\mathbb{A}$-valued $S$-family $\left(
a_{s}\right)  _{s\in S}$ satisfy $\left\vert S\right\vert =0$. If you are
using Definition \ref{def.ind.gen-com.defsum2} instead of Definition
\ref{def.ind.gen-com.defsum1}, you should instead be using Lemma
\ref{lem.ind.gen-com.same-sum} \textbf{(a)} to argue this.)

Comparing (\ref{pf.prop.ind.gen-com.n(n+1)/2.IB.1}) with $\dfrac{0\left(
0+1\right)  }{2}=0$, we obtain $\sum_{i\in\left\{  1,2,\ldots,0\right\}
}i=\dfrac{0\left(  0+1\right)  }{2}$. In other words,
(\ref{eq.prop.ind.gen-com.n(n+1)/2.claim}) holds for $n=0$. This completes the
induction base.

\textit{Induction step:} Let $m\in\mathbb{N}$. Assume that
(\ref{eq.prop.ind.gen-com.n(n+1)/2.claim}) holds for $n=m$. We must prove that
(\ref{eq.prop.ind.gen-com.n(n+1)/2.claim}) holds for $n=m+1$.

We have assumed that (\ref{eq.prop.ind.gen-com.n(n+1)/2.claim}) holds for
$n=m$. In other words, we have
\begin{equation}
\sum_{i\in\left\{  1,2,\ldots,m\right\}  }i=\dfrac{m\left(  m+1\right)  }{2}.
\label{pf.prop.ind.gen-com.n(n+1)/2.IH}%
\end{equation}

Now, $\left\vert \left\{  1,2,\ldots,m+1\right\}  \right\vert =m+1$ and
$m+1\in\left\{  1,2,\ldots,m+1\right\}  $ (since $m+1$ is a positive integer
(since $m\in\mathbb{N}$)). Hence, (\ref{eq.sum.def.1}) (applied to $n=m$,
$S=\left\{  1,2,\ldots,m+1\right\}  $, $t=m+1$ and $\left(  a_{s}\right)
_{s\in S}=\left(  i\right)  _{i\in\left\{  1,2,\ldots,m+1\right\}  }$) yields%
\begin{equation}
\sum_{i\in\left\{  1,2,\ldots,m+1\right\}  }i=\left(  m+1\right)  +\sum
_{i\in\left\{  1,2,\ldots,m+1\right\}  \setminus\left\{  m+1\right\}  }i.
\label{pf.prop.ind.gen-com.n(n+1)/2.3}%
\end{equation}
(Here, we have relied on the equality (\ref{eq.sum.def.1}), which appears
verbatim in Definition \ref{def.ind.gen-com.defsum1}. If you are using
Definition \ref{def.ind.gen-com.defsum2} instead of Definition
\ref{def.ind.gen-com.defsum1}, you should instead be using Lemma
\ref{lem.ind.gen-com.same-sum} \textbf{(b)} to argue this.)

Now, (\ref{pf.prop.ind.gen-com.n(n+1)/2.3}) becomes%
\begin{align*}
\sum_{i\in\left\{  1,2,\ldots,m+1\right\}  }i  &  =\left(  m+1\right)
+\sum_{i\in\left\{  1,2,\ldots,m+1\right\}  \setminus\left\{  m+1\right\}
}i=\left(  m+1\right)  +\underbrace{\sum_{i\in\left\{  1,2,\ldots,m\right\}
}i}_{\substack{=\dfrac{m\left(  m+1\right)  }{2}\\\text{(by
(\ref{pf.prop.ind.gen-com.n(n+1)/2.IH}))}}}\\
&  \ \ \ \ \ \ \ \ \ \ \left(  \text{since }\left\{  1,2,\ldots,m+1\right\}
\setminus\left\{  m+1\right\}  =\left\{  1,2,\ldots,m\right\}  \right) \\
&  =\left(  m+1\right)  +\dfrac{m\left(  m+1\right)  }{2}=\dfrac{2\left(
m+1\right)  +m\left(  m+1\right)  }{2}\\
&  =\dfrac{\left(  m+1\right)  \left(  \left(  m+1\right)  +1\right)  }{2}%
\end{align*}
(since $2\left(  m+1\right)  +m\left(  m+1\right)  =\left(  m+1\right)
\left(  \left(  m+1\right)  +1\right)  $). In other words,
(\ref{eq.prop.ind.gen-com.n(n+1)/2.claim}) holds for $n=m+1$. This completes
the induction step. Thus, the induction proof of
(\ref{eq.prop.ind.gen-com.n(n+1)/2.claim}) is finished. Hence, Proposition
\ref{prop.ind.gen-com.n(n+1)/2} holds.
\end{proof}

\subsubsection{Sums of a few numbers}

Merely for the sake of future convenience, let us restate (\ref{eq.sum.def.1})
in a slightly more direct way (without mentioning $\left\vert S\right\vert $):

\begin{proposition}
\label{prop.ind.gen-com.split-off}Let $S$ be a finite set, and let $\left(
a_{s}\right)  _{s\in S}$ be an $\mathbb{A}$-valued $S$-family. Let $t\in S$.
Then,%
\[
\sum_{s\in S}a_{s}=a_{t}+\sum_{s\in S\setminus\left\{  t\right\}  }a_{s}.
\]

\end{proposition}

\begin{proof}
[Proof of Proposition \ref{prop.ind.gen-com.split-off}.]Let $n=\left\vert
S\setminus\left\{  t\right\}  \right\vert $; thus, $n\in\mathbb{N}$ (since
$S\setminus\left\{  t\right\}  $ is a finite set). Also, $n=\left\vert
S\setminus\left\{  t\right\}  \right\vert =\left\vert S\right\vert -1$ (since
$t\in S$), and thus $\left\vert S\right\vert =n+1$. Hence, (\ref{eq.sum.def.1}%
) yields $\sum_{s\in S}a_{s}=a_{t}+\sum_{s\in S\setminus\left\{  t\right\}
}a_{s}$. This proves Proposition \ref{prop.ind.gen-com.split-off}.
\end{proof}

(Alternatively, we can argue that Proposition \ref{prop.ind.gen-com.split-off}
is the same as Lemma \ref{lem.ind.gen-com.same-sum} \textbf{(b)}, except that
we are now using Definition \ref{def.ind.gen-com.defsum1} instead of
Definition \ref{def.ind.gen-com.defsum2} to define the sums involved -- but
this difference is insubstantial, since we have shown that these two
definitions are equivalent.)

\begin{noncompile}
(Proof of Proposition \ref{prop.ind.gen-com.split-off}.) If we use Definition
\ref{def.ind.gen-com.defsum2} instead of Definition
\ref{def.ind.gen-com.defsum1} (for defining sums), then the claim of
Proposition \ref{prop.ind.gen-com.split-off} becomes precisely the claim of
Lemma \ref{lem.ind.gen-com.same-sum} \textbf{(b)}. Since we know that
Definition \ref{def.ind.gen-com.defsum2} is equivalent to Definition
\ref{def.ind.gen-com.defsum1} (according to Theorem \ref{thm.ind.gen-com.wd}
\textbf{(b)}), we thus conclude that Proposition
\ref{prop.ind.gen-com.split-off} is equivalent to the claim of Lemma
\ref{lem.ind.gen-com.same-sum} \textbf{(b)}. Thus, Proposition
\ref{prop.ind.gen-com.split-off} holds (since Lemma
\ref{lem.ind.gen-com.same-sum} \textbf{(b)} holds).
\end{noncompile}

In Section \ref{sect.sums-repetitorium}, we have introduced $a_{u}%
+a_{u+1}+\cdots+a_{v}$ as an abbreviation for the sum $\sum_{s=u}^{v}%
a_{s}=\sum_{s\in\left\{  u,u+1,\ldots,v\right\}  }a_{s}$ (whenever $u$ and $v$
are two integers, and $a_{s}$ is a number for each $s\in\left\{
u,u+1,\ldots,v\right\}  $). In order to ensure that this abbreviation does not
create any nasty surprises, we need to check that it behaves as we would
expect -- i.e., that it satisfies the following four properties:

\begin{itemize}
\item If the sum $a_{u}+a_{u+1}+\cdots+a_{v}$ has no addends (i.e., if $u>v$),
then it equals $0$.

\item If the sum $a_{u}+a_{u+1}+\cdots+a_{v}$ has exactly one addend (i.e., if
$u=v$), then it equals $a_{u}$.

\item If the sum $a_{u}+a_{u+1}+\cdots+a_{v}$ has exactly two addends (i.e.,
if $u=v-1$), then it equals $a_{u}+a_{v}$.

\item If $v\geq u$, then
\begin{align*}
a_{u}+a_{u+1}+\cdots+a_{v}  &  =\left(  a_{u}+a_{u+1}+\cdots+a_{v-1}\right)
+a_{v}\\
&  =a_{u}+\left(  a_{u+1}+a_{u+2}+\cdots+a_{v}\right)  .
\end{align*}

\end{itemize}

The first of these four properties follows from the definition (indeed, if
$u>v$, then the set $\left\{  u,u+1,\ldots,v\right\}  $ is empty and thus
satisfies $\left\vert \left\{  u,u+1,\ldots,v\right\}  \right\vert =0$; but
this yields $\sum_{s\in\left\{  u,u+1,\ldots,v\right\}  }a_{s}=0$). The fourth
of these four properties can easily be obtained from Proposition
\ref{prop.ind.gen-com.split-off}\footnote{In more detail: Assume that $v\geq
u$. Thus, both $u$ and $v$ belong to the set $\left\{  u,u+1,\ldots,v\right\}
$. Hence, Proposition \ref{prop.ind.gen-com.split-off} (applied to $S=\left\{
u,u+1,\ldots,v\right\}  $ and $t=v$) yields $\sum_{s\in\left\{  u,u+1,\ldots
,v\right\}  }a_{s}=a_{v}+\sum_{s\in\left\{  u,u+1,\ldots,v\right\}
\setminus\left\{  v\right\}  }a_{s}$. Thus,%
\begin{align*}
&  a_{u}+a_{u+1}+\cdots+a_{v}\\
&  =\sum_{s\in\left\{  u,u+1,\ldots,v\right\}  }a_{s}=a_{v}+\sum_{s\in\left\{
u,u+1,\ldots,v\right\}  \setminus\left\{  v\right\}  }a_{s}\\
&  =a_{v}+\sum_{s\in\left\{  u,u+1,\ldots,v-1\right\}  }a_{s}%
\ \ \ \ \ \ \ \ \ \ \left(  \text{since }\left\{  u,u+1,\ldots,v\right\}
\setminus\left\{  v\right\}  =\left\{  u,u+1,\ldots,v-1\right\}  \right) \\
&  =a_{v}+\left(  a_{u}+a_{u+1}+\cdots+a_{v-1}\right)  =\left(  a_{u}%
+a_{u+1}+\cdots+a_{v-1}\right)  +a_{v}.
\end{align*}
\par
Also, Proposition \ref{prop.ind.gen-com.split-off} (applied to $S=\left\{
u,u+1,\ldots,v\right\}  $ and $t=u$) yields $\sum_{s\in\left\{  u,u+1,\ldots
,v\right\}  }a_{s}=a_{u}+\sum_{s\in\left\{  u,u+1,\ldots,v\right\}
\setminus\left\{  u\right\}  }a_{s}$. Thus,%
\begin{align*}
&  a_{u}+a_{u+1}+\cdots+a_{v}\\
&  =\sum_{s\in\left\{  u,u+1,\ldots,v\right\}  }a_{s}=a_{u}+\sum_{s\in\left\{
u,u+1,\ldots,v\right\}  \setminus\left\{  u\right\}  }a_{s}\\
&  =a_{u}+\sum_{s\in\left\{  u+1,u+2,\ldots,v\right\}  }a_{s}%
\ \ \ \ \ \ \ \ \ \ \left(  \text{since }\left\{  u,u+1,\ldots,v\right\}
\setminus\left\{  u\right\}  =\left\{  u+1,u+2,\ldots,v\right\}  \right) \\
&  =a_{u}+\left(  a_{u+1}+a_{u+2}+\cdots+a_{v}\right)  .
\end{align*}
Hence,%
\[
a_{u}+a_{u+1}+\cdots+a_{v}=\left(  a_{u}+a_{u+1}+\cdots+a_{v-1}\right)
+a_{v}=a_{u}+\left(  a_{u+1}+a_{u+2}+\cdots+a_{v}\right)  .
\]
}. The second and third properties follow from the following fact:

\begin{proposition}
\label{prop.ind.gen-com.sum12}Let $S$ be a finite set. For every $s\in S$, let
$a_{s}$ be an element of $\mathbb{A}$.

\textbf{(a)} If $S=\left\{  p\right\}  $ for some element $p$, then this $p$
satisfies%
\[
\sum_{s\in S}a_{s}=a_{p}.
\]

\textbf{(b)} If $S=\left\{  p,q\right\}  $ for two distinct elements $p$ and
$q$, then these $p$ and $q$ satisfy
\[
\sum_{s\in S}a_{s}=a_{p}+a_{q}.
\]

\end{proposition}

\begin{proof}
[Proof of Proposition \ref{prop.ind.gen-com.sum12}.]\textbf{(a)} Assume that
$S=\left\{  p\right\}  $ for some element $p$. Consider this $p$.

The first bullet point of Definition \ref{def.ind.gen-com.defsum1} shows that
$\sum_{s\in\varnothing}a_{s}=0$ (since $\left\vert \varnothing\right\vert
=0$). But $p\in\left\{  p\right\}  =S$. Hence, Proposition
\ref{prop.ind.gen-com.split-off} (applied to $t=p$) yields
\begin{align*}
\sum_{s\in S}a_{s}  &  =a_{p}+\sum_{s\in S\setminus\left\{  p\right\}  }%
a_{s}=a_{p}+\underbrace{\sum_{s\in\varnothing}a_{s}}_{=0}%
\ \ \ \ \ \ \ \ \ \ \left(  \text{since }\underbrace{S}_{=\left\{  p\right\}
}\setminus\left\{  p\right\}  =\left\{  p\right\}  \setminus\left\{
p\right\}  =\varnothing\right) \\
&  =a_{p}+0=a_{p}.
\end{align*}
This proves Proposition \ref{prop.ind.gen-com.sum12} \textbf{(a)}.

\textbf{(b)} Assume that $S=\left\{  p,q\right\}  $ for two distinct elements
$p$ and $q$. Consider these $p$ and $q$. Thus, $q\neq p$ (since $p$ and $q$
are distinct), so that $q\notin\left\{  p\right\}  $.

Proposition \ref{prop.ind.gen-com.sum12} \textbf{(a)} (applied to $\left\{
p\right\}  $ instead of $S$) yields $\sum_{s\in\left\{  p\right\}  }%
a_{s}=a_{p}$ (since $\left\{  p\right\}  =\left\{  p\right\}  $).

We have $\underbrace{S}_{=\left\{  p,q\right\}  =\left\{  p\right\}
\cup\left\{  q\right\}  }\setminus\left\{  q\right\}  =\left(  \left\{
p\right\}  \cup\left\{  q\right\}  \right)  \setminus\left\{  q\right\}
=\left\{  p\right\}  \setminus\left\{  q\right\}  =\left\{  p\right\}  $
(since $q\notin\left\{  p\right\}  $). Also, $q\in\left\{  p,q\right\}  =S$.
Hence, Proposition \ref{prop.ind.gen-com.split-off} (applied to $t=q$) yields
\begin{align*}
\sum_{s\in S}a_{s}  &  =a_{q}+\sum_{s\in S\setminus\left\{  q\right\}  }%
a_{s}=a_{q}+\underbrace{\sum_{s\in\left\{  p\right\}  }a_{s}}_{=a_{p}\text{ }%
}\ \ \ \ \ \ \ \ \ \ \left(  \text{since }S\setminus\left\{  q\right\}
=\left\{  p\right\}  \right) \\
&  =a_{q}+a_{p}=a_{p}+a_{q}.
\end{align*}
This proves Proposition \ref{prop.ind.gen-com.sum12} \textbf{(b)}.
\end{proof}

\subsubsection{Linearity of sums}

We shall now prove some general properties of finite sums. We begin with the
equality (\ref{eq.sum.linear1}) from Section \ref{sect.sums-repetitorium}:

\begin{theorem}
\label{thm.ind.gen-com.sum(a+b)}Let $S$ be a finite set. For every $s\in S$,
let $a_{s}$ and $b_{s}$ be elements of $\mathbb{A}$. Then,%
\[
\sum_{s\in S}\left(  a_{s}+b_{s}\right)  =\sum_{s\in S}a_{s}+\sum_{s\in
S}b_{s}.
\]

\end{theorem}

Before we prove this theorem, let us show a simple lemma:

\begin{lemma}
\label{lem.ind.gen-com.xyuv}Let $x$, $y$, $u$ and $v$ be four numbers (i.e.,
elements of $\mathbb{A}$). Then,%
\[
\left(  x+y\right)  +\left(  u+v\right)  =\left(  x+u\right)  +\left(
y+v\right)  .
\]

\end{lemma}

\begin{proof}
[Proof of Lemma \ref{lem.ind.gen-com.xyuv}.]Proposition
\ref{prop.ind.gen-com.fgh} (applied to $a=y$, $b=u$ and $c=v$) yields%
\begin{equation}
\left(  y+u\right)  +v=y+\left(  u+v\right)  .
\label{pf.lem.ind.gen-com.xyuv.1}%
\end{equation}
Also, Proposition \ref{prop.ind.gen-com.fgh} (applied to $a=x$, $b=y$ and
$c=u+v$) yields%
\begin{align}
\left(  x+y\right)  +\left(  u+v\right)   &  =x+\underbrace{\left(  y+\left(
u+v\right)  \right)  }_{\substack{=\left(  y+u\right)  +v\\\text{(by
(\ref{pf.lem.ind.gen-com.xyuv.1}))}}}\nonumber\\
&  =x+\left(  \left(  y+u\right)  +v\right)  .
\label{pf.lem.ind.gen-com.xyuv.2}%
\end{align}
The same argument (with $y$ and $u$ replaced by $u$ and $y$) yields%
\begin{equation}
\left(  x+u\right)  +\left(  y+v\right)  =x+\left(  \left(  u+y\right)
+v\right)  . \label{pf.lem.ind.gen-com.xyuv.3}%
\end{equation}
But Proposition \ref{prop.ind.gen-com.fg} (applied to $a=y$ and $b=u$) yields
$y+u=u+y$. Thus, (\ref{pf.lem.ind.gen-com.xyuv.2}) becomes%
\[
\left(  x+y\right)  +\left(  u+v\right)  =x+\left(  \underbrace{\left(
y+u\right)  }_{=u+y}+v\right)  =x+\left(  \left(  u+y\right)  +v\right)
=\left(  x+u\right)  +\left(  y+v\right)
\]
(by (\ref{pf.lem.ind.gen-com.xyuv.3})). This proves Lemma
\ref{lem.ind.gen-com.xyuv}.
\end{proof}

\begin{proof}
[Proof of Theorem \ref{thm.ind.gen-com.sum(a+b)}.]Forget that we fixed $S$,
$a_{s}$ and $b_{s}$. We shall prove Theorem \ref{thm.ind.gen-com.sum(a+b)} by
induction on $\left\vert S\right\vert $:

\textit{Induction base:} Theorem \ref{thm.ind.gen-com.sum(a+b)} holds under
the condition that $\left\vert S\right\vert =0$%
\ \ \ \ \footnote{\textit{Proof.} Let $S$, $a_{s}$ and $b_{s}$ be as in
Theorem \ref{thm.ind.gen-com.sum(a+b)}. Assume that $\left\vert S\right\vert
=0$. Thus, the first bullet point of Definition \ref{def.ind.gen-com.defsum1}
yields $\sum_{s\in S}a_{s}=0$ and $\sum_{s\in S}b_{s}=0$ and $\sum_{s\in
S}\left(  a_{s}+b_{s}\right)  =0$. Hence,%
\[
\sum_{s\in S}\left(  a_{s}+b_{s}\right)  =0=\underbrace{0}_{=\sum_{s\in
S}a_{s}}+\underbrace{0}_{=\sum_{s\in S}b_{s}}=\sum_{s\in S}a_{s}+\sum_{s\in
S}b_{s}.
\]
\par
Now, forget that we fixed $S$, $a_{s}$ and $b_{s}$. We thus have proved that
if $S$, $a_{s}$ and $b_{s}$ are as in Theorem \ref{thm.ind.gen-com.sum(a+b)},
and if $\left\vert S\right\vert =0$, then $\sum_{s\in S}\left(  a_{s}%
+b_{s}\right)  =\sum_{s\in S}a_{s}+\sum_{s\in S}b_{s}$. In other words,
Theorem \ref{thm.ind.gen-com.sum(a+b)} holds under the condition that
$\left\vert S\right\vert =0$. Qed.}. This completes the induction base.

\textit{Induction step:} Let $m\in\mathbb{N}$. Assume that Theorem
\ref{thm.ind.gen-com.sum(a+b)} holds under the condition that $\left\vert
S\right\vert =m$. We must now prove that Theorem
\ref{thm.ind.gen-com.sum(a+b)} holds under the condition that $\left\vert
S\right\vert =m+1$.

We have assumed that Theorem \ref{thm.ind.gen-com.sum(a+b)} holds under the
condition that $\left\vert S\right\vert =m$. In other words, the following
claim holds:

\begin{statement}
\textit{Claim 1:} Let $S$ be a finite set such that $\left\vert S\right\vert
=m$. For every $s\in S$, let $a_{s}$ and $b_{s}$ be elements of $\mathbb{A}$.
Then,%
\[
\sum_{s\in S}\left(  a_{s}+b_{s}\right)  =\sum_{s\in S}a_{s}+\sum_{s\in
S}b_{s}.
\]

\end{statement}

Next, we shall show the following claim:

\begin{statement}
\textit{Claim 2:} Let $S$ be a finite set such that $\left\vert S\right\vert
=m+1$. For every $s\in S$, let $a_{s}$ and $b_{s}$ be elements of $\mathbb{A}%
$. Then,%
\[
\sum_{s\in S}\left(  a_{s}+b_{s}\right)  =\sum_{s\in S}a_{s}+\sum_{s\in
S}b_{s}.
\]

\end{statement}

[\textit{Proof of Claim 2:} We have $\left\vert S\right\vert =m+1>m\geq0$.
Hence, the set $S$ is nonempty. Thus, there exists some $t\in S$. Consider
this $t$.

From $t\in S$, we obtain $\left\vert S\setminus\left\{  t\right\}  \right\vert
=\left\vert S\right\vert -1=m$ (since $\left\vert S\right\vert =m+1$). Hence,
Claim 1 (applied to $S\setminus\left\{  t\right\}  $ instead of $S$) yields%
\begin{equation}
\sum_{s\in S\setminus\left\{  t\right\}  }\left(  a_{s}+b_{s}\right)
=\sum_{s\in S\setminus\left\{  t\right\}  }a_{s}+\sum_{s\in S\setminus\left\{
t\right\}  }b_{s}. \label{pf.thm.ind.gen-com.sum(a+b).c2.pf.1}%
\end{equation}

Now, Proposition \ref{prop.ind.gen-com.split-off} yields%
\begin{equation}
\sum_{s\in S}a_{s}=a_{t}+\sum_{s\in S\setminus\left\{  t\right\}  }a_{s}.
\label{pf.thm.ind.gen-com.sum(a+b).c2.pf.a}%
\end{equation}
Also, Proposition \ref{prop.ind.gen-com.split-off} (applied to $b_{s}$ instead
of $a_{s}$) yields%
\begin{equation}
\sum_{s\in S}b_{s}=b_{t}+\sum_{s\in S\setminus\left\{  t\right\}  }b_{s}.
\label{pf.thm.ind.gen-com.sum(a+b).c2.pf.b}%
\end{equation}
Finally, Proposition \ref{prop.ind.gen-com.split-off} (applied to $a_{s}%
+b_{s}$ instead of $a_{s}$) yields%
\begin{align*}
\sum_{s\in S}\left(  a_{s}+b_{s}\right)   &  =\left(  a_{t}+b_{t}\right)
+\underbrace{\sum_{s\in S\setminus\left\{  t\right\}  }\left(  a_{s}%
+b_{s}\right)  }_{\substack{=\sum_{s\in S\setminus\left\{  t\right\}  }%
a_{s}+\sum_{s\in S\setminus\left\{  t\right\}  }b_{s}\\\text{(by
(\ref{pf.thm.ind.gen-com.sum(a+b).c2.pf.1}))}}}\\
&  =\left(  a_{t}+b_{t}\right)  +\left(  \sum_{s\in S\setminus\left\{
t\right\}  }a_{s}+\sum_{s\in S\setminus\left\{  t\right\}  }b_{s}\right) \\
&  =\underbrace{\left(  a_{t}+\sum_{s\in S\setminus\left\{  t\right\}  }%
a_{s}\right)  }_{\substack{=\sum_{s\in S}a_{s}\\\text{(by
(\ref{pf.thm.ind.gen-com.sum(a+b).c2.pf.a}))}}}+\underbrace{\left(  b_{t}%
+\sum_{s\in S\setminus\left\{  t\right\}  }b_{s}\right)  }_{\substack{=\sum
_{s\in S}b_{s}\\\text{(by (\ref{pf.thm.ind.gen-com.sum(a+b).c2.pf.b}))}}}\\
&  \ \ \ \ \ \ \ \ \ \ \left(
\begin{array}
[c]{c}%
\text{by Lemma \ref{lem.ind.gen-com.xyuv} (applied}\\
\text{to }x=a_{t}\text{, }y=b_{t}\text{, }u=\sum_{s\in S\setminus\left\{
t\right\}  }a_{s}\text{ and }v=\sum_{s\in S\setminus\left\{  t\right\}  }%
b_{s}\text{)}%
\end{array}
\right) \\
&  =\sum_{s\in S}a_{s}+\sum_{s\in S}b_{s}.
\end{align*}
This proves Claim 2.]

But Claim 2 says precisely that Theorem \ref{thm.ind.gen-com.sum(a+b)} holds
under the condition that $\left\vert S\right\vert =m+1$. Hence, we conclude
that Theorem \ref{thm.ind.gen-com.sum(a+b)} holds under the condition that
$\left\vert S\right\vert =m+1$ (since Claim 2 is proven). This completes the
induction step. Thus, Theorem \ref{thm.ind.gen-com.sum(a+b)} is proven by induction.
\end{proof}

We shall next prove (\ref{eq.sum.linear2}):

\begin{theorem}
\label{thm.ind.gen-com.sum(la)}Let $S$ be a finite set. For every $s\in S$,
let $a_{s}$ be an element of $\mathbb{A}$. Also, let $\lambda$ be an element
of $\mathbb{A}$. Then,%
\[
\sum_{s\in S}\lambda a_{s}=\lambda\sum_{s\in S}a_{s}.
\]

\end{theorem}

To prove this theorem, we need the following fundamental fact of arithmetic:

\begin{proposition}
\label{prop.ind.gen-com.distr}Let $x$, $y$ and $z$ be three numbers (i.e.,
elements of $\mathbb{A}$). Then, $x\left(  y+z\right)  =xy+xz$.
\end{proposition}

Proposition \ref{prop.ind.gen-com.distr} is known as the
\textit{distributivity} (or \textit{left distributivity}) in $\mathbb{A}$. It
is a fundamental result, and its proof can be found in standard
textbooks\footnote{For example, Proposition \ref{prop.ind.gen-com.distr} is
proven in \cite[Theorem 3.2.3 (6)]{Swanso18} for the case when $\mathbb{A}%
=\mathbb{N}$; in \cite[Theorem 3.5.4 (6)]{Swanso18} for the case when
$\mathbb{A}=\mathbb{Z}$; in \cite[Theorem 3.6.4 (6)]{Swanso18} for the case
when $\mathbb{A}=\mathbb{Q}$; in \cite[Theorem 3.7.14]{Swanso18} for the case
when $\mathbb{A}=\mathbb{R}$; in \cite[Theorem 3.9.2]{Swanso18} for the case
when $\mathbb{A}=\mathbb{C}$.}.

\begin{proof}
[Proof of Theorem \ref{thm.ind.gen-com.sum(la)}.]Forget that we fixed $S$,
$a_{s}$ and $\lambda$. We shall prove Theorem \ref{thm.ind.gen-com.sum(la)} by
induction on $\left\vert S\right\vert $:

\begin{vershort}
\textit{Induction base:} The induction base (i.e., proving that Theorem
\ref{thm.ind.gen-com.sum(la)} holds under the condition that $\left\vert
S\right\vert =0$) is similar to the induction base in the proof of Theorem
\ref{thm.ind.gen-com.sum(a+b)} above; we thus leave it to the reader.
\end{vershort}

\begin{verlong}
\textit{Induction base:} Theorem \ref{thm.ind.gen-com.sum(la)} holds under the
condition that $\left\vert S\right\vert =0$\ \ \ \ \footnote{\textit{Proof.}
Let $S$, $a_{s}$ and $\lambda$ be as in Theorem \ref{thm.ind.gen-com.sum(la)}.
Assume that $\left\vert S\right\vert =0$. Thus, the first bullet point of
Definition \ref{def.ind.gen-com.defsum1} yields $\sum_{s\in S}a_{s}=0$ and
$\sum_{s\in S}\lambda a_{s}=0$. Hence,%
\[
\sum_{s\in S}\lambda a_{s}=0=\lambda\underbrace{0}_{=\sum_{s\in S}a_{s}%
}=\lambda\sum_{s\in S}a_{s}.
\]
\par
Now, forget that we fixed $S$, $a_{s}$ and $\lambda$. We thus have proved that
if $S$, $a_{s}$ and $\lambda$ are as in Theorem \ref{thm.ind.gen-com.sum(la)},
and if $\left\vert S\right\vert =0$, then $\sum_{s\in S}\lambda a_{s}%
=\lambda\sum_{s\in S}a_{s}$. In other words, Theorem
\ref{thm.ind.gen-com.sum(la)} holds under the condition that $\left\vert
S\right\vert =0$. Qed.}. This completes the induction base.
\end{verlong}

\textit{Induction step:} Let $m\in\mathbb{N}$. Assume that Theorem
\ref{thm.ind.gen-com.sum(la)} holds under the condition that $\left\vert
S\right\vert =m$. We must now prove that Theorem \ref{thm.ind.gen-com.sum(la)}
holds under the condition that $\left\vert S\right\vert =m+1$.

We have assumed that Theorem \ref{thm.ind.gen-com.sum(la)} holds under the
condition that $\left\vert S\right\vert =m$. In other words, the following
claim holds:

\begin{statement}
\textit{Claim 1:} Let $S$ be a finite set such that $\left\vert S\right\vert
=m$. For every $s\in S$, let $a_{s}$ be an element of $\mathbb{A}$. Also, let
$\lambda$ be an element of $\mathbb{A}$. Then,%
\[
\sum_{s\in S}\lambda a_{s}=\lambda\sum_{s\in S}a_{s}.
\]

\end{statement}

Next, we shall show the following claim:

\begin{statement}
\textit{Claim 2:} Let $S$ be a finite set such that $\left\vert S\right\vert
=m+1$. For every $s\in S$, let $a_{s}$ be an element of $\mathbb{A}$. Also,
let $\lambda$ be an element of $\mathbb{A}$. Then,%
\[
\sum_{s\in S}\lambda a_{s}=\lambda\sum_{s\in S}a_{s}.
\]

\end{statement}

[\textit{Proof of Claim 2:} We have $\left\vert S\right\vert =m+1>m\geq0$.
Hence, the set $S$ is nonempty. Thus, there exists some $t\in S$. Consider
this $t$.

From $t\in S$, we obtain $\left\vert S\setminus\left\{  t\right\}  \right\vert
=\left\vert S\right\vert -1=m$ (since $\left\vert S\right\vert =m+1$). Hence,
Claim 1 (applied to $S\setminus\left\{  t\right\}  $ instead of $S$) yields%
\begin{equation}
\sum_{s\in S\setminus\left\{  t\right\}  }\lambda a_{s}=\lambda\sum_{s\in
S\setminus\left\{  t\right\}  }a_{s}.
\label{pf.thm.ind.gen-com.sum(la).c2.pf.1}%
\end{equation}

Now, Proposition \ref{prop.ind.gen-com.split-off} yields%
\[
\sum_{s\in S}a_{s}=a_{t}+\sum_{s\in S\setminus\left\{  t\right\}  }a_{s}.
\]
Multiplying both sides of this equality by $\lambda$, we obtain%
\[
\lambda\sum_{s\in S}a_{s}=\lambda\left(  a_{t}+\sum_{s\in S\setminus\left\{
t\right\}  }a_{s}\right)  =\lambda a_{t}+\lambda\sum_{s\in S\setminus\left\{
t\right\}  }a_{s}%
\]
(by Proposition \ref{prop.ind.gen-com.distr} (applied to $x=\lambda$,
$y=a_{t}$ and $z=\sum_{s\in S\setminus\left\{  t\right\}  }a_{s}$)). Also,
Proposition \ref{prop.ind.gen-com.split-off} (applied to $\lambda a_{s}$
instead of $a_{s}$) yields%
\[
\sum_{s\in S}\lambda a_{s}=\lambda a_{t}+\underbrace{\sum_{s\in S\setminus
\left\{  t\right\}  }\lambda a_{s}}_{\substack{=\lambda\sum_{s\in
S\setminus\left\{  t\right\}  }a_{s}\\\text{(by
(\ref{pf.thm.ind.gen-com.sum(la).c2.pf.1}))}}}=\lambda a_{t}+\lambda\sum_{s\in
S\setminus\left\{  t\right\}  }a_{s}.
\]
Comparing the preceding two equalities, we find%
\[
\sum_{s\in S}\lambda a_{s}=\lambda\sum_{s\in S}a_{s}.
\]
This proves Claim 2.]

But Claim 2 says precisely that Theorem \ref{thm.ind.gen-com.sum(la)} holds
under the condition that $\left\vert S\right\vert =m+1$. Hence, we conclude
that Theorem \ref{thm.ind.gen-com.sum(la)} holds under the condition that
$\left\vert S\right\vert =m+1$ (since Claim 2 is proven). This completes the
induction step. Thus, Theorem \ref{thm.ind.gen-com.sum(la)} is proven by induction.
\end{proof}

Finally, let us prove (\ref{eq.sum.sum0}):

\begin{theorem}
\label{thm.ind.gen-com.sum(0)}Let $S$ be a finite set. Then,
\[
\sum_{s\in S}0=0.
\]

\end{theorem}

\begin{proof}
[Proof of Theorem \ref{thm.ind.gen-com.sum(0)}.]It is completely
straightforward to prove Theorem \ref{thm.ind.gen-com.sum(0)} by induction on
$\left\vert S\right\vert $ (as we proved Theorem \ref{thm.ind.gen-com.sum(la)}%
, for example). But let us give an even shorter argument: Theorem
\ref{thm.ind.gen-com.sum(la)} (applied to $a_{s}=0$ and $\lambda=0$) yields%
\[
\sum_{s\in S}0\cdot0=0\sum_{s\in S}0=0.
\]
In view of $0\cdot0=0$, this rewrites as $\sum_{s\in S}0=0$. This proves
Theorem \ref{thm.ind.gen-com.sum(0)}.
\end{proof}

\subsubsection{Splitting a sum by a value of a function}

We shall now prove a more complicated (but crucial) property of finite sums --
namely, the equality (\ref{eq.sum.sheph}) in the case when $W$ is
finite\footnote{We prefer to only treat the case when $W$ is finite for now.
The case when $W$ is infinite would require us to properly introduce the
notion of an infinite sum with only finitely many nonzero terms. While this is
not hard to do, we aren't quite ready for it yet (see Theorem
\ref{thm.ind.gen-com.sheph} further below for this).}:

\begin{theorem}
\label{thm.ind.gen-com.shephf}Let $S$ be a finite set. Let $W$ be a finite
set. Let $f:S\rightarrow W$ be a map. Let $a_{s}$ be an element of
$\mathbb{A}$ for each $s\in S$. Then,%
\[
\sum_{s\in S}a_{s}=\sum_{w\in W}\sum_{\substack{s\in S;\\f\left(  s\right)
=w}}a_{s}.
\]

\end{theorem}

Here, we are using the following convention (made in Section
\ref{sect.sums-repetitorium}):

\begin{convention}
Let $S$ be a finite set. Let $\mathcal{A}\left(  s\right)  $ be a logical
statement defined for every $s\in S$. For each $s\in S$ satisfying
$\mathcal{A}\left(  s\right)  $, let $a_{s}$ be a number (i.e., an element of
$\mathbb{A}$). Then, we set%
\[
\sum_{\substack{s\in S;\\\mathcal{A}\left(  s\right)  }}a_{s}=\sum
_{s\in\left\{  t\in S\ \mid\ \mathcal{A}\left(  t\right)  \right\}  }a_{s}.
\]

\end{convention}

Thus, the sum $\sum_{\substack{s\in S;\\f\left(  s\right)  =w}}a_{s}$ in
Theorem \ref{thm.ind.gen-com.shephf} can be rewritten as $\sum_{s\in\left\{
t\in S\ \mid\ f\left(  t\right)  =w\right\}  }a_{s}$.

Our proof of Theorem \ref{thm.ind.gen-com.shephf} relies on the following
simple set-theoretic fact:

\begin{lemma}
\label{lem.ind.gen-com.shephf-l1}Let $S$ and $W$ be two sets. Let
$f:S\rightarrow W$ be a map. Let $q\in S$. Let $g$ be the restriction
$f\mid_{S\setminus\left\{  q\right\}  }$ of the map $f$ to $S\setminus\left\{
q\right\}  $. Let $w\in W$. Then,%
\[
\left\{  t\in S\setminus\left\{  q\right\}  \ \mid\ g\left(  t\right)
=w\right\}  =\left\{  t\in S\ \mid\ f\left(  t\right)  =w\right\}
\setminus\left\{  q\right\}  .
\]

\end{lemma}

\begin{vershort}
\begin{proof}
[Proof of Lemma \ref{lem.ind.gen-com.shephf-l1}.]We know that $g$ is the
restriction $f\mid_{S\setminus\left\{  q\right\}  }$ of the map $f$ to
$S\setminus\left\{  q\right\}  $. Thus, $g$ is a map from $S\setminus\left\{
q\right\}  $ to $W$ and satisfies%
\begin{equation}
g\left(  t\right)  =f\left(  t\right)  \ \ \ \ \ \ \ \ \ \ \text{for each
}t\in S\setminus\left\{  q\right\}  .
\label{pf.lem.ind.gen-com.shephf-l1.short.1}%
\end{equation}

Now,%
\begin{align*}
&  \left\{  t\in S\setminus\left\{  q\right\}  \ \mid\ \underbrace{g\left(
t\right)  }_{\substack{=f\left(  t\right)  \\\text{(by
(\ref{pf.lem.ind.gen-com.shephf-l1.short.1}))}}}=w\right\} \\
&  =\left\{  t\in S\setminus\left\{  q\right\}  \ \mid\ f\left(  t\right)
=w\right\}  =\left\{  t\in S\ \mid\ f\left(  t\right)  =w\text{ and }t\in
S\setminus\left\{  q\right\}  \right\} \\
&  =\left\{  t\in S\ \mid\ f\left(  t\right)  =w\right\}  \cap
\underbrace{\left\{  t\in S\ \mid\ t\in S\setminus\left\{  q\right\}
\right\}  }_{=S\setminus\left\{  q\right\}  }\\
&  =\left\{  t\in S\ \mid\ f\left(  t\right)  =w\right\}  \cap\left(
S\setminus\left\{  q\right\}  \right)  =\left\{  t\in S\ \mid\ f\left(
t\right)  =w\right\}  \setminus\left\{  q\right\}  .
\end{align*}
This proves Lemma \ref{lem.ind.gen-com.shephf-l1}.
\end{proof}
\end{vershort}

\begin{verlong}
\begin{proof}
[Proof of Lemma \ref{lem.ind.gen-com.shephf-l1}.]We know that $g$ is the
restriction $f\mid_{S\setminus\left\{  q\right\}  }$ of the map $f$ to
$S\setminus\left\{  q\right\}  $. Thus, $g$ is a map from $S\setminus\left\{
q\right\}  $ to $W$ and satisfies%
\begin{equation}
\left(  g\left(  s\right)  =f\left(  s\right)  \ \ \ \ \ \ \ \ \ \ \text{for
each }s\in S\setminus\left\{  q\right\}  \right)  .
\label{pf.lem.ind.gen-com.shephf-l1.long.1}%
\end{equation}

Let $p\in\left\{  t\in S\setminus\left\{  q\right\}  \ \mid\ g\left(
t\right)  =w\right\}  $. Thus, $p$ is a $t\in S\setminus\left\{  q\right\}  $
satisfying $g\left(  t\right)  =w$. In other words, $p$ is an element of
$S\setminus\left\{  q\right\}  $ and satisfies $g\left(  p\right)  =w$.
Applying (\ref{pf.lem.ind.gen-com.shephf-l1.long.1}) to $s=p$, we obtain
$g\left(  p\right)  =f\left(  p\right)  $. Hence, $f\left(  p\right)
=g\left(  p\right)  =w$. Thus, $p$ is an element of $S$ (since $p\in
S\setminus\left\{  q\right\}  \subseteq S$) and satisfies $f\left(  p\right)
=w$. In other words, $p$ is a $t\in S$ satisfying $f\left(  t\right)  =w$. In
other words, $p\in\left\{  t\in S\ \mid\ f\left(  t\right)  =w\right\}  $.
Moreover, $p\notin\left\{  q\right\}  $ (since $p\in S\setminus\left\{
q\right\}  $). Combining $p\in\left\{  t\in S\ \mid\ f\left(  t\right)
=w\right\}  $ with $p\notin\left\{  q\right\}  $, we obtain $p\in\left\{  t\in
S\ \mid\ f\left(  t\right)  =w\right\}  \setminus\left\{  q\right\}  $.

Now, forget that we fixed $p$. We thus have proven that $p\in\left\{  t\in
S\ \mid\ f\left(  t\right)  =w\right\}  \setminus\left\{  q\right\}  $ for
each $p\in\left\{  t\in S\setminus\left\{  q\right\}  \ \mid\ g\left(
t\right)  =w\right\}  $. In other words,%
\begin{equation}
\left\{  t\in S\setminus\left\{  q\right\}  \ \mid\ g\left(  t\right)
=w\right\}  \subseteq\left\{  t\in S\ \mid\ f\left(  t\right)  =w\right\}
\setminus\left\{  q\right\}  . \label{pf.lem.ind.gen-com.shephf-l1.long.4}%
\end{equation}

On the other hand, let $s\in\left\{  t\in S\ \mid\ f\left(  t\right)
=w\right\}  \setminus\left\{  q\right\}  $. Thus, $s\in\left\{  t\in
S\ \mid\ f\left(  t\right)  =w\right\}  $ and $s\notin\left\{  q\right\}  $.
In particular, we have $s\in\left\{  t\in S\ \mid\ f\left(  t\right)
=w\right\}  $. In other words, $s$ is a $t\in S$ satisfying $f\left(
t\right)  =w$. In other words, $s$ is an element of $S$ and satisfies
$f\left(  s\right)  =w$. Since $s$ is an element of $S$, we have $s\in S$.
Combining this with $s\notin\left\{  q\right\}  $, we obtain $s\in
S\setminus\left\{  q\right\}  $. Thus,
(\ref{pf.lem.ind.gen-com.shephf-l1.long.1}) yields $g\left(  s\right)
=f\left(  s\right)  =w$. Hence, $s$ is an element of $S\setminus\left\{
q\right\}  $ (since $s\in S\setminus\left\{  q\right\}  $) and satisfies
$g\left(  s\right)  =w$. In other words, $s$ is a $t\in S\setminus\left\{
q\right\}  $ satisfying $g\left(  t\right)  =w$. In other words, $s\in\left\{
t\in S\setminus\left\{  q\right\}  \ \mid\ g\left(  t\right)  =w\right\}  $.

Now, forget that we fixed $s$. We thus have shown that $s\in\left\{  t\in
S\setminus\left\{  q\right\}  \ \mid\ g\left(  t\right)  =w\right\}  $ for
each $s\in\left\{  t\in S\ \mid\ f\left(  t\right)  =w\right\}  \setminus
\left\{  q\right\}  $. In other words,%
\[
\left\{  t\in S\ \mid\ f\left(  t\right)  =w\right\}  \setminus\left\{
q\right\}  \subseteq\left\{  t\in S\setminus\left\{  q\right\}  \ \mid
\ g\left(  t\right)  =w\right\}  .
\]
Combining this with (\ref{pf.lem.ind.gen-com.shephf-l1.long.4}), we obtain%
\begin{equation}
\left\{  t\in S\setminus\left\{  q\right\}  \ \mid\ g\left(  t\right)
=w\right\}  =\left\{  t\in S\ \mid\ f\left(  t\right)  =w\right\}
\setminus\left\{  q\right\}  .\nonumber
\end{equation}
This proves Lemma \ref{lem.ind.gen-com.shephf-l1}.
\end{proof}
\end{verlong}

\begin{proof}
[Proof of Theorem \ref{thm.ind.gen-com.shephf}.]We shall prove Theorem
\ref{thm.ind.gen-com.shephf} by induction on $\left\vert S\right\vert $:

\textit{Induction base:} Theorem \ref{thm.ind.gen-com.shephf} holds under the
condition that $\left\vert S\right\vert =0$\ \ \ \ \footnote{\textit{Proof.}
Let $S$, $W$, $f$ and $a_{s}$ be as in Theorem \ref{thm.ind.gen-com.shephf}.
Assume that $\left\vert S\right\vert =0$. Thus, the first bullet point of
Definition \ref{def.ind.gen-com.defsum1} yields $\sum_{s\in S}a_{s}=0$.
Moreover, $S=\varnothing$ (since $\left\vert S\right\vert =0$). Hence, each
$w\in W$ satisfies%
\begin{align*}
\sum_{\substack{s\in S;\\f\left(  s\right)  =w}}a_{s}  &  =\sum_{s\in\left\{
t\in S\ \mid\ f\left(  t\right)  =w\right\}  }a_{s}=\sum_{s\in\varnothing
}a_{s}\\
&  \ \ \ \ \ \ \ \ \ \ \left(
\begin{array}
[c]{c}%
\text{since }\left\{  t\in S\ \mid\ f\left(  s\right)  =w\right\}
=\varnothing\\
\text{(because }\left\{  t\in S\ \mid\ f\left(  s\right)  =w\right\}
\subseteq S=\varnothing\text{)}%
\end{array}
\right) \\
&  =\left(  \text{empty sum}\right)  =0.
\end{align*}
Summing these equalities over all $w\in W$, we obtain%
\[
\sum_{w\in W}\sum_{\substack{s\in S;\\f\left(  s\right)  =w}}a_{s}=\sum_{w\in
W}0=0
\]
(by an application of Theorem \ref{thm.ind.gen-com.sum(0)}). Comparing this
with $\sum_{s\in S}a_{s}=0$, we obtain $\sum_{s\in S}a_{s}=\sum_{w\in W}%
\sum_{\substack{s\in S;\\f\left(  s\right)  =w}}a_{s}$.
\par
Now, forget that we fixed $S$, $W$, $f$ and $a_{s}$. We thus have proved that
if $S$, $W$, $f$ and $a_{s}$ are as in Theorem \ref{thm.ind.gen-com.shephf},
and if $\left\vert S\right\vert =0$, then $\sum_{s\in S}a_{s}=\sum_{w\in
W}\sum_{\substack{s\in S;\\f\left(  s\right)  =w}}a_{s}$. In other words,
Theorem \ref{thm.ind.gen-com.shephf} holds under the condition that
$\left\vert S\right\vert =0$. Qed.}. This completes the induction base.

\textit{Induction step:} Let $m\in\mathbb{N}$. Assume that Theorem
\ref{thm.ind.gen-com.shephf} holds under the condition that $\left\vert
S\right\vert =m$. We must now prove that Theorem \ref{thm.ind.gen-com.shephf}
holds under the condition that $\left\vert S\right\vert =m+1$.

We have assumed that Theorem \ref{thm.ind.gen-com.shephf} holds under the
condition that $\left\vert S\right\vert =m$. In other words, the following
claim holds:

\begin{statement}
\textit{Claim 1:} Let $S$ be a finite set such that $\left\vert S\right\vert
=m$. Let $W$ be a finite set. Let $f:S\rightarrow W$ be a map. Let $a_{s}$ be
an element of $\mathbb{A}$ for each $s\in S$. Then,%
\[
\sum_{s\in S}a_{s}=\sum_{w\in W}\sum_{\substack{s\in S;\\f\left(  s\right)
=w}}a_{s}.
\]

\end{statement}

Next, we shall show the following claim:

\begin{statement}
\textit{Claim 2:} Let $S$ be a finite set such that $\left\vert S\right\vert
=m+1$. Let $W$ be a finite set. Let $f:S\rightarrow W$ be a map. Let $a_{s}$
be an element of $\mathbb{A}$ for each $s\in S$. Then,%
\[
\sum_{s\in S}a_{s}=\sum_{w\in W}\sum_{\substack{s\in S;\\f\left(  s\right)
=w}}a_{s}.
\]

\end{statement}

[\textit{Proof of Claim 2:} We have $\left\vert S\right\vert =m+1>m\geq0$.
Hence, the set $S$ is nonempty. Thus, there exists some $q\in S$. Consider
this $q$.

From $q\in S$, we obtain $\left\vert S\setminus\left\{  q\right\}  \right\vert
=\left\vert S\right\vert -1=m$ (since $\left\vert S\right\vert =m+1$).

Let $g$ be the restriction $f\mid_{S\setminus\left\{  q\right\}  }$ of the map
$f$ to $S\setminus\left\{  q\right\}  $. Thus, $g$ is a map from
$S\setminus\left\{  q\right\}  $ to $W$.

For each $w\in W$, we define a number $b_{w}$ by%
\begin{equation}
b_{w}=\sum_{\substack{s\in S;\\f\left(  s\right)  =w}}a_{s}.
\label{pf.lem.ind.gen-com.shephf.c2.bw=}%
\end{equation}

Furthermore, for each $w\in W$, we define a number $c_{w}$ by%
\begin{equation}
c_{w}=\sum_{\substack{s\in S\setminus\left\{  q\right\}  ;\\g\left(  s\right)
=w}}a_{s}. \label{pf.lem.ind.gen-com.shephf.c2.cw=}%
\end{equation}

Recall that $\left\vert S\setminus\left\{  q\right\}  \right\vert =m$. Hence,
Claim 1 (applied to $S\setminus\left\{  q\right\}  $ and $g$ instead of $S$
and $f$) yields%
\begin{equation}
\sum_{s\in S\setminus\left\{  q\right\}  }a_{s}=\sum_{w\in W}\underbrace{\sum
_{\substack{s\in S\setminus\left\{  q\right\}  ;\\g\left(  s\right)  =w}%
}a_{s}}_{\substack{=c_{w}\\\text{(by (\ref{pf.lem.ind.gen-com.shephf.c2.cw=}%
))}}}=\sum_{w\in W}c_{w}. \label{pf.lem.ind.gen-com.shephf.c2.usec1}%
\end{equation}

Every $w\in W\setminus\left\{  f\left(  q\right)  \right\}  $ satisfies%
\begin{equation}
b_{w}=c_{w}. \label{pf.lem.ind.gen-com.shephf.c2.4b}%
\end{equation}

[\textit{Proof of (\ref{pf.lem.ind.gen-com.shephf.c2.4b}):} Let $w\in
W\setminus\left\{  f\left(  q\right)  \right\}  $. Thus, $w\in W$ and
$w\notin\left\{  f\left(  q\right)  \right\}  $.

\begin{vershort}
If we had $q\in\left\{  t\in S\ \mid\ f\left(  t\right)  =w\right\}  $, then
we would have $f\left(  q\right)  =w$, which would lead to $w=f\left(
q\right)  \in\left\{  f\left(  q\right)  \right\}  $; but this would
contradict $w\notin\left\{  f\left(  q\right)  \right\}  $. Hence, we cannot
have $q\in\left\{  t\in S\ \mid\ f\left(  t\right)  =w\right\}  $. Hence, we
have $q\notin\left\{  t\in S\ \mid\ f\left(  t\right)  =w\right\}  $.
\end{vershort}

\begin{verlong}
Let us next prove that $q\notin\left\{  t\in S\ \mid\ f\left(  t\right)
=w\right\}  $. Indeed, assume the contrary (for the sake of contradiction).
Thus, $q\in\left\{  t\in S\ \mid\ f\left(  t\right)  =w\right\}  $. In other
words, $q$ is a $t\in S$ satisfying $f\left(  t\right)  =w$. In other words,
$q$ is an element of $S$ and satisfies $f\left(  q\right)  =w$. Hence,
$w=f\left(  q\right)  \in\left\{  f\left(  q\right)  \right\}  $; but this
contradicts $w\notin\left\{  f\left(  q\right)  \right\}  $.

This contradiction shows that our assumption was false. Hence, $q\notin%
\left\{  t\in S\ \mid\ f\left(  t\right)  =w\right\}  $ is proven.
\end{verlong}

But $w\in W$; thus, Lemma \ref{lem.ind.gen-com.shephf-l1} yields%
\begin{align}
\left\{  t\in S\setminus\left\{  q\right\}  \ \mid\ g\left(  t\right)
=w\right\}   &  =\left\{  t\in S\ \mid\ f\left(  t\right)  =w\right\}
\setminus\left\{  q\right\} \nonumber\\
&  =\left\{  t\in S\ \mid\ f\left(  t\right)  =w\right\}
\label{pf.lem.ind.gen-com.shephf.c2.4b.pf.2}%
\end{align}
(since $q\notin\left\{  t\in S\ \mid\ f\left(  t\right)  =w\right\}  $).

On the other hand, the definition of $b_{w}$ yields%
\begin{equation}
b_{w}=\sum_{\substack{s\in S;\\f\left(  s\right)  =w}}a_{s}=\sum_{s\in\left\{
t\in S\ \mid\ f\left(  t\right)  =w\right\}  }a_{s}
\label{pf.lem.ind.gen-com.shephf.c2.4b.pf.1}%
\end{equation}
(by the definition of the \textquotedblleft$\sum_{\substack{s\in S;\\f\left(
s\right)  =w}}$\textquotedblright\ symbol). Also, the definition of $c_{w}$
yields%
\begin{align*}
c_{w}  &  =\sum_{\substack{s\in S\setminus\left\{  q\right\}  ;\\g\left(
s\right)  =w}}a_{s}=\sum_{s\in\left\{  t\in S\setminus\left\{  q\right\}
\ \mid\ g\left(  t\right)  =w\right\}  }a_{s}=\sum_{s\in\left\{  t\in
S\ \mid\ f\left(  t\right)  =w\right\}  }a_{s}\\
&  \ \ \ \ \ \ \ \ \ \ \left(
\begin{array}
[c]{c}%
\text{since }\left\{  t\in S\setminus\left\{  q\right\}  \ \mid\ g\left(
t\right)  =w\right\}  =\left\{  t\in S\ \mid\ f\left(  t\right)  =w\right\} \\
\text{(by (\ref{pf.lem.ind.gen-com.shephf.c2.4b.pf.2}))}%
\end{array}
\right) \\
&  =b_{w}\ \ \ \ \ \ \ \ \ \ \left(  \text{by
(\ref{pf.lem.ind.gen-com.shephf.c2.4b.pf.1})}\right)  .
\end{align*}
Thus, $b_{w}=c_{w}$. This proves (\ref{pf.lem.ind.gen-com.shephf.c2.4b}).]

Also,
\begin{equation}
b_{f\left(  q\right)  }=a_{q}+c_{f\left(  q\right)  }.
\label{pf.lem.ind.gen-com.shephf.c2.5b}%
\end{equation}

[\textit{Proof of (\ref{pf.lem.ind.gen-com.shephf.c2.5b}):} Define a subset
$U$ of $S$ by%
\begin{equation}
U=\left\{  t\in S\ \mid\ f\left(  t\right)  =f\left(  q\right)  \right\}  .
\label{pf.lem.ind.gen-com.shephf.c2.5b.pf.U=}%
\end{equation}

\begin{verlong}
This set $U$ is a subset of $S$, and thus is finite (since $S$ is finite).
\end{verlong}

We can apply Lemma \ref{lem.ind.gen-com.shephf-l1} to $w=f\left(  q\right)  $.
We thus obtain%
\begin{align}
\left\{  t\in S\setminus\left\{  q\right\}  \ \mid\ g\left(  t\right)
=f\left(  q\right)  \right\}   &  =\underbrace{\left\{  t\in S\ \mid\ f\left(
t\right)  =f\left(  q\right)  \right\}  }_{\substack{=U\\\text{(by
(\ref{pf.lem.ind.gen-com.shephf.c2.5b.pf.U=}))}}}\setminus\left\{  q\right\}
\nonumber\\
&  =U\setminus\left\{  q\right\}  .
\label{pf.lem.ind.gen-com.shephf.c2.5b.pf.1}%
\end{align}

We know that $q$ is a $t\in S$ satisfying $f\left(  t\right)  =f\left(
q\right)  $ (since $q\in S$ and $f\left(  q\right)  =f\left(  q\right)  $). In
other words, $q\in\left\{  t\in S\ \mid\ f\left(  t\right)  =f\left(
q\right)  \right\}  $. In other words, $q\in U$ (since $U=\left\{  t\in
S\ \mid\ f\left(  t\right)  =f\left(  q\right)  \right\}  $). Thus,
Proposition \ref{prop.ind.gen-com.split-off} (applied to $U$ and $q$ instead
of $S$ and $t$) yields%
\begin{equation}
\sum_{s\in U}a_{s}=a_{q}+\sum_{s\in U\setminus\left\{  q\right\}  }a_{s}.
\label{pf.lem.ind.gen-com.shephf.c2.5b.pf.2}%
\end{equation}

But (\ref{pf.lem.ind.gen-com.shephf.c2.5b.pf.1}) shows that $U\setminus
\left\{  q\right\}  =\left\{  t\in S\setminus\left\{  q\right\}
\ \mid\ g\left(  t\right)  =f\left(  q\right)  \right\}  $. Thus,%
\begin{equation}
\sum_{s\in U\setminus\left\{  q\right\}  }a_{s}=\sum_{s\in\left\{  t\in
S\setminus\left\{  q\right\}  \ \mid\ g\left(  t\right)  =f\left(  q\right)
\right\}  }a_{s}=c_{f\left(  q\right)  }
\label{pf.lem.ind.gen-com.shephf.c2.5b.pf.3}%
\end{equation}
(since the definition of $c_{f\left(  q\right)  }$ yields $c_{f\left(
q\right)  }=\sum_{\substack{s\in S\setminus\left\{  q\right\}  ;\\g\left(
s\right)  =f\left(  q\right)  }}a_{s}=\sum_{s\in\left\{  t\in S\setminus
\left\{  q\right\}  \ \mid\ g\left(  t\right)  =f\left(  q\right)  \right\}
}a_{s}$).

On the other hand, the definition of $b_{f\left(  q\right)  }$ yields%
\begin{align*}
b_{f\left(  q\right)  }  &  =\sum_{\substack{s\in S;\\f\left(  s\right)
=f\left(  q\right)  }}a_{s}=\sum_{s\in\left\{  t\in S\ \mid\ f\left(
t\right)  =f\left(  q\right)  \right\}  }a_{s}\\
&  =\sum_{s\in U}a_{s}\ \ \ \ \ \ \ \ \ \ \left(  \text{since }\left\{  t\in
S\ \mid\ f\left(  t\right)  =f\left(  q\right)  \right\}  =U\right) \\
&  =a_{q}+\underbrace{\sum_{s\in U\setminus\left\{  q\right\}  }a_{s}%
}_{\substack{=c_{f\left(  q\right)  }\\\text{(by
(\ref{pf.lem.ind.gen-com.shephf.c2.5b.pf.3}))}}}\ \ \ \ \ \ \ \ \ \ \left(
\text{by (\ref{pf.lem.ind.gen-com.shephf.c2.5b.pf.2})}\right) \\
&  =a_{q}+c_{f\left(  q\right)  }.
\end{align*}
This proves (\ref{pf.lem.ind.gen-com.shephf.c2.5b}).]

Now, recall that $q\in S$. Hence, Proposition \ref{prop.ind.gen-com.split-off}
(applied to $t=q$) yields%
\begin{equation}
\sum_{s\in S}a_{s}=a_{q}+\sum_{s\in S\setminus\left\{  q\right\}  }a_{s}.
\label{pf.lem.ind.gen-com.shephf.c2.6}%
\end{equation}

Also, $f\left(  q\right)  \in W$. Hence, Proposition
\ref{prop.ind.gen-com.split-off} (applied to $W$, $\left(  c_{w}\right)
_{w\in W}$ and $f\left(  q\right)  $ instead of $S$, $\left(  a_{s}\right)
_{s\in S}$ and $t$) yields%
\[
\sum_{w\in W}c_{w}=c_{f\left(  q\right)  }+\sum_{w\in W\setminus\left\{
f\left(  q\right)  \right\}  }c_{w}.
\]
Hence, (\ref{pf.lem.ind.gen-com.shephf.c2.usec1}) becomes%
\begin{equation}
\sum_{s\in S\setminus\left\{  q\right\}  }a_{s}=\sum_{w\in W}c_{w}=c_{f\left(
q\right)  }+\sum_{w\in W\setminus\left\{  f\left(  q\right)  \right\}  }c_{w}.
\label{pf.lem.ind.gen-com.shephf.c2.8}%
\end{equation}

Also, Proposition \ref{prop.ind.gen-com.split-off} (applied to $W$, $\left(
b_{w}\right)  _{w\in W}$ and $f\left(  q\right)  $ instead of $S$, $\left(
a_{s}\right)  _{s\in S}$ and $t$) yields%
\begin{align*}
\sum_{w\in W}b_{w}  &  =\underbrace{b_{f\left(  q\right)  }}_{\substack{=a_{q}%
+c_{f\left(  q\right)  }\\\text{(by (\ref{pf.lem.ind.gen-com.shephf.c2.5b}))}%
}}+\sum_{w\in W\setminus\left\{  f\left(  q\right)  \right\}  }%
\underbrace{b_{w}}_{\substack{=c_{w}\\\text{(by
(\ref{pf.lem.ind.gen-com.shephf.c2.4b}))}}}\\
&  =\left(  a_{q}+c_{f\left(  q\right)  }\right)  +\sum_{w\in W\setminus
\left\{  f\left(  q\right)  \right\}  }c_{w}=a_{q}+\left(  c_{f\left(
q\right)  }+\sum_{w\in W\setminus\left\{  f\left(  q\right)  \right\}  }%
c_{w}\right)
\end{align*}
(by Proposition \ref{prop.ind.gen-com.fgh}, applied to $a_{q}$, $c_{f\left(
q\right)  }$ and $\sum_{w\in W\setminus\left\{  f\left(  q\right)  \right\}
}c_{w}$ instead of $a$, $b$ and $c$). Thus,%
\[
\sum_{w\in W}b_{w}=a_{q}+\underbrace{\left(  c_{f\left(  q\right)  }%
+\sum_{w\in W\setminus\left\{  f\left(  q\right)  \right\}  }c_{w}\right)
}_{\substack{=\sum_{s\in S\setminus\left\{  q\right\}  }a_{s}\\\text{(by
(\ref{pf.lem.ind.gen-com.shephf.c2.8}))}}}=a_{q}+\sum_{s\in S\setminus\left\{
q\right\}  }a_{s}=\sum_{s\in S}a_{s}%
\]
(by (\ref{pf.lem.ind.gen-com.shephf.c2.6})). Hence,%
\[
\sum_{s\in S}a_{s}=\sum_{w\in W}\underbrace{b_{w}}_{\substack{=\sum
_{\substack{s\in S;\\f\left(  s\right)  =w}}a_{s}\\\text{(by
(\ref{pf.lem.ind.gen-com.shephf.c2.bw=}))}}}=\sum_{w\in W}\sum_{\substack{s\in
S;\\f\left(  s\right)  =w}}a_{s}.
\]
This proves Claim 2.]

But Claim 2 says precisely that Theorem \ref{thm.ind.gen-com.shephf} holds
under the condition that $\left\vert S\right\vert =m+1$. Hence, we conclude
that Theorem \ref{thm.ind.gen-com.shephf} holds under the condition that
$\left\vert S\right\vert =m+1$ (since Claim 2 is proven). This completes the
induction step. Thus, Theorem \ref{thm.ind.gen-com.shephf} is proven by induction.
\end{proof}

\subsubsection{Splitting a sum into two}

Next, we shall prove the equality (\ref{eq.sum.split}):

\begin{theorem}
\label{thm.ind.gen-com.split2}Let $S$ be a finite set. Let $X$ and $Y$ be two
subsets of $S$ such that $X\cap Y=\varnothing$ and $X\cup Y=S$. (Equivalently,
$X$ and $Y$ are two subsets of $S$ such that each element of $S$ lies in
\textbf{exactly} one of $X$ and $Y$.) Let $a_{s}$ be a number (i.e., an
element of $\mathbb{A}$) for each $s\in S$. Then,%
\[
\sum_{s\in S}a_{s}=\sum_{s\in X}a_{s}+\sum_{s\in Y}a_{s}.
\]

\end{theorem}

\begin{proof}
[Proof of Theorem \ref{thm.ind.gen-com.split2}.]From the assumptions $X\cap
Y=\varnothing$ and $X\cup Y=S$, we can easily obtain $S\setminus X=Y$.

\begin{verlong}
[\textit{Proof:} The set $X\cap Y$ is empty (since $X\cap Y=\varnothing$);
thus, it has no elements.

Let $y\in Y$. If we had $y\in X$, then we would have $y\in X\cap Y$ (since
$y\in X$ and $y\in Y$), which would show that the set $X\cap Y$ has at least
one element (namely, $y$); but this would contradict the fact that this set
$X\cap Y$ has no elements. Thus, we cannot have $y\in X$. Hence, we have
$y\notin X$. Combining $y\in Y\subseteq S$ with $y\notin X$, we obtain $y\in
S\setminus X$.

Now, forget that we fixed $y$. We thus have shown that $y\in S\setminus X$ for
each $y\in Y$. In other words, $Y\subseteq S\setminus X$.

Combining this with%
\[
\underbrace{S}_{=X\cup Y}\setminus X=\left(  X\cup Y\right)  \setminus
X=Y\setminus X\subseteq Y,
\]
we obtain $S\setminus X=Y$. Qed.]
\end{verlong}

We define a map $f:S\rightarrow\left\{  0,1\right\}  $ by setting%
\[
\left(  f\left(  s\right)  =%
\begin{cases}
0, & \text{if }s\in X;\\
1, & \text{if }s\notin X
\end{cases}
\ \ \ \ \ \ \ \ \ \ \text{for every }s\in S\right)  .
\]

For each $w\in\left\{  0,1\right\}  $, we define a number $b_{w}$ by%
\begin{equation}
b_{w}=\sum_{\substack{s\in S;\\f\left(  s\right)  =w}}a_{s}.
\label{pf.thm.ind.gen-com.split2.bw=}%
\end{equation}

\begin{vershort}
Proposition \ref{prop.ind.gen-com.sum12} \textbf{(b)} (applied to $\left\{
0,1\right\}  $, $0$, $1$ and $\left(  b_{w}\right)  _{w\in\left\{
0,1\right\}  }$ instead of $S$, $p$, $q$ and $\left(  a_{s}\right)  _{s\in S}%
$) yields $\sum_{w\in\left\{  0,1\right\}  }b_{w}=b_{0}+b_{1}$.
\end{vershort}

\begin{verlong}
We have $\left\{  0,1\right\}  =\left\{  0,1\right\}  $, where $0$ and $1$ are
two distinct elements. Hence, Proposition \ref{prop.ind.gen-com.sum12}
\textbf{(b)} (applied to $\left\{  0,1\right\}  $, $0$, $1$ and $\left(
b_{w}\right)  _{w\in\left\{  0,1\right\}  }$ instead of $S$, $p$, $q$ and
$\left(  a_{s}\right)  _{s\in S}$) yields $\sum_{w\in\left\{  0,1\right\}
}b_{w}=b_{0}+b_{1}$.
\end{verlong}

Now, Theorem \ref{thm.ind.gen-com.shephf} (applied to $W=\left\{  0,1\right\}
$) yields
\begin{equation}
\sum_{s\in S}a_{s}=\sum_{w\in\left\{  0,1\right\}  }\underbrace{\sum
_{\substack{s\in S;\\f\left(  s\right)  =w}}a_{s}}_{\substack{=b_{w}%
\\\text{(by (\ref{pf.thm.ind.gen-com.split2.bw=}))}}}=\sum_{w\in\left\{
0,1\right\}  }b_{w}=b_{0}+b_{1}. \label{pf.thm.ind.gen-com.split2.1}%
\end{equation}

On the other hand,
\begin{equation}
b_{0}=\sum_{s\in X}a_{s}. \label{pf.thm.ind.gen-com.split2.b0=}%
\end{equation}

\begin{vershort}
[\textit{Proof of (\ref{pf.thm.ind.gen-com.split2.b0=}):} The definition of
the map $f$ shows that an element $t\in S$ satisfies $f\left(  t\right)  =0$
\textbf{if and only if} it belongs to $X$. Hence, the set of all elements
$t\in S$ that satisfy $f\left(  t\right)  =0$ is precisely $X$. In other
words,%
\[
\left\{  t\in S\ \mid\ f\left(  t\right)  =0\right\}  =X.
\]
But the definition of $b_{0}$ yields%
\[
b_{0}=\sum_{\substack{s\in S;\\f\left(  s\right)  =0}}a_{s}=\sum_{s\in\left\{
t\in S\ \mid\ f\left(  t\right)  =0\right\}  }a_{s}=\sum_{s\in X}a_{s}%
\]
(since $\left\{  t\in S\ \mid\ f\left(  t\right)  =0\right\}  =X$). This
proves (\ref{pf.thm.ind.gen-com.split2.b0=}).]
\end{vershort}

\begin{verlong}
[\textit{Proof of (\ref{pf.thm.ind.gen-com.split2.b0=}):} Let $p\in X$. Thus,
the definition of $f$ yields $f\left(  p\right)  =%
\begin{cases}
0, & \text{if }p\in X;\\
1, & \text{if }p\notin X
\end{cases}
=0$ (since $p\in X$). Hence, $p$ is a $t\in S$ satisfying $f\left(  t\right)
=0$ (since $p\in S$ and $f\left(  p\right)  =0$). In other words,
$p\in\left\{  t\in S\ \mid\ f\left(  t\right)  =0\right\}  $.

Now, forget that we fixed $p$. We thus have proven that $p\in\left\{  t\in
S\ \mid\ f\left(  t\right)  =0\right\}  $ for each $p\in X$. In other words,
\begin{equation}
X\subseteq\left\{  t\in S\ \mid\ f\left(  t\right)  =0\right\}  .
\label{pf.thm.ind.gen-com.split2.b0=.pf.1}%
\end{equation}

On the other hand, let $s\in\left\{  t\in S\ \mid\ f\left(  t\right)
=0\right\}  $. Thus, $s$ is a $t\in S$ satisfying $f\left(  t\right)  =0$. In
other words, $s$ is an element of $S$ and satisfies $f\left(  s\right)  =0$.
If we had $s\notin X$, then we would have%
\begin{align*}
f\left(  s\right)   &  =%
\begin{cases}
0, & \text{if }s\in X;\\
1, & \text{if }s\notin X
\end{cases}
\ \ \ \ \ \ \ \ \ \ \left(  \text{by the definition of }f\right) \\
&  =1\ \ \ \ \ \ \ \ \ \ \left(  \text{since }s\notin X\right)  ,
\end{align*}
which would contradict $f\left(  s\right)  =0\neq1$. Hence, we cannot have
$s\notin X$. Thus, we have $s\in X$.

Now, forget that we fixed $s$. We thus have proven that $s\in X$ for each
$s\in\left\{  t\in S\ \mid\ f\left(  t\right)  =0\right\}  $. In other words,
$\left\{  t\in S\ \mid\ f\left(  t\right)  =0\right\}  \subseteq X$. Combining
this with (\ref{pf.thm.ind.gen-com.split2.b0=.pf.1}), we obtain%
\[
X=\left\{  t\in S\ \mid\ f\left(  t\right)  =0\right\}  .
\]
Hence,%
\[
\sum_{s\in X}a_{s}=\sum_{s\in\left\{  t\in S\ \mid\ f\left(  t\right)
=0\right\}  }a_{s}.
\]
Comparing this with%
\begin{align*}
b_{0}  &  =\sum_{\substack{s\in S;\\f\left(  s\right)  =0}}a_{s}%
\ \ \ \ \ \ \ \ \ \ \left(  \text{by the definition of }b_{0}\right) \\
&  =\sum_{s\in\left\{  t\in S\ \mid\ f\left(  t\right)  =0\right\}  }a_{s},
\end{align*}
we obtain $b_{0}=\sum_{s\in X}a_{s}$. This proves
(\ref{pf.thm.ind.gen-com.split2.b0=}).]
\end{verlong}

Furthermore,
\begin{equation}
b_{1}=\sum_{s\in Y}a_{s}. \label{pf.thm.ind.gen-com.split2.b1=}%
\end{equation}

\begin{vershort}
[\textit{Proof of (\ref{pf.thm.ind.gen-com.split2.b1=}):} The definition of
the map $f$ shows that an element $t\in S$ satisfies $f\left(  t\right)  =1$
\textbf{if and only if} $t\notin X$. Thus, for each $t\in S$, we have the
following chain of equivalences:%
\[
\left(  f\left(  t\right)  =1\right)  \ \Longleftrightarrow\ \left(  t\notin
X\right)  \ \Longleftrightarrow\ \left(  t\in S\setminus X\right)
\ \Longleftrightarrow\ \left(  t\in Y\right)
\]
(since $S\setminus X=Y$). In other words, an element $t\in S$ satisfies
$f\left(  t\right)  =1$ \textbf{if and only if} $t$ belongs to $Y$. Hence, the
set of all elements $t\in S$ that satisfy $f\left(  t\right)  =1$ is precisely
$Y$. In other words,%
\[
\left\{  t\in S\ \mid\ f\left(  t\right)  =1\right\}  =Y.
\]
But the definition of $b_{1}$ yields%
\[
b_{1}=\sum_{\substack{s\in S;\\f\left(  s\right)  =1}}a_{s}=\sum_{s\in\left\{
t\in S\ \mid\ f\left(  t\right)  =1\right\}  }a_{s}=\sum_{s\in Y}a_{s}%
\]
(since $\left\{  t\in S\ \mid\ f\left(  t\right)  =1\right\}  =Y$). This
proves (\ref{pf.thm.ind.gen-com.split2.b1=}).]
\end{vershort}

\begin{verlong}
[\textit{Proof of (\ref{pf.thm.ind.gen-com.split2.b1=}):} Let $p\in Y$. Thus,
$p\in Y=S\setminus X$ (since $S\setminus X=Y$). In other words, $p\in S$ and
$p\notin X$. The definition of $f$ yields $f\left(  p\right)  =%
\begin{cases}
0, & \text{if }p\in X;\\
1, & \text{if }p\notin X
\end{cases}
=1$ (since $p\notin X$). Hence, $p$ is a $t\in S$ satisfying $f\left(
t\right)  =1$ (since $p\in S$ and $f\left(  p\right)  =1$). In other words,
$p\in\left\{  t\in S\ \mid\ f\left(  t\right)  =1\right\}  $.

Now, forget that we fixed $p$. We thus have proven that $p\in\left\{  t\in
S\ \mid\ f\left(  t\right)  =1\right\}  $ for each $p\in Y$. In other words,
\begin{equation}
Y\subseteq\left\{  t\in S\ \mid\ f\left(  t\right)  =1\right\}  .
\label{pf.thm.ind.gen-com.split2.b1=.pf.1}%
\end{equation}

On the other hand, let $s\in\left\{  t\in S\ \mid\ f\left(  t\right)
=1\right\}  $. Thus, $s$ is a $t\in S$ satisfying $f\left(  t\right)  =1$. In
other words, $s$ is an element of $S$ and satisfies $f\left(  s\right)  =1$.
If we had $s\in X$, then we would have%
\begin{align*}
f\left(  s\right)   &  =%
\begin{cases}
0, & \text{if }s\in X;\\
1, & \text{if }s\notin X
\end{cases}
\ \ \ \ \ \ \ \ \ \ \left(  \text{by the definition of }f\right) \\
&  =0\ \ \ \ \ \ \ \ \ \ \left(  \text{since }s\in X\right)  ,
\end{align*}
which would contradict $f\left(  s\right)  =1\neq0$. Hence, we cannot have
$s\in X$. Thus, we have $s\notin X$. Combining $s\in\left\{  t\in
S\ \mid\ f\left(  t\right)  =1\right\}  \subseteq S$ with $s\notin X$, we
obtain $s\in S\setminus X=Y$.

Now, forget that we fixed $s$. We thus have proven that $s\in Y$ for each
$s\in\left\{  t\in S\ \mid\ f\left(  t\right)  =1\right\}  $. In other words,
$\left\{  t\in S\ \mid\ f\left(  t\right)  =1\right\}  \subseteq Y$. Combining
this with (\ref{pf.thm.ind.gen-com.split2.b1=.pf.1}), we obtain%
\[
Y=\left\{  t\in S\ \mid\ f\left(  t\right)  =1\right\}  .
\]
Hence,%
\[
\sum_{s\in Y}a_{s}=\sum_{s\in\left\{  t\in S\ \mid\ f\left(  t\right)
=1\right\}  }a_{s}.
\]
Comparing this with%
\begin{align*}
b_{1}  &  =\sum_{\substack{s\in S;\\f\left(  s\right)  =1}}a_{s}%
\ \ \ \ \ \ \ \ \ \ \left(  \text{by the definition of }b_{1}\right) \\
&  =\sum_{s\in\left\{  t\in S\ \mid\ f\left(  t\right)  =1\right\}  }a_{s},
\end{align*}
we obtain $b_{1}=\sum_{s\in Y}a_{s}$. This proves
(\ref{pf.thm.ind.gen-com.split2.b1=}).]
\end{verlong}

Now, (\ref{pf.thm.ind.gen-com.split2.1}) becomes%
\[
\sum_{s\in S}a_{s}=\underbrace{b_{0}}_{\substack{=\sum_{s\in X}a_{s}%
\\\text{(by (\ref{pf.thm.ind.gen-com.split2.b0=}))}}}+\underbrace{b_{1}%
}_{\substack{=\sum_{s\in Y}a_{s}\\\text{(by
(\ref{pf.thm.ind.gen-com.split2.b1=}))}}}=\sum_{s\in X}a_{s}+\sum_{s\in
Y}a_{s}.
\]
This proves Theorem \ref{thm.ind.gen-com.split2}.
\end{proof}

Similarly, we can prove the equality (\ref{eq.sum.split-n}). (This proof was
already outlined in Section \ref{sect.sums-repetitorium}.)

A consequence of Theorem \ref{thm.ind.gen-com.split2} is the following fact,
which has appeared as the equality (\ref{eq.sum.drop0}) in Section
\ref{sect.sums-repetitorium}:

\begin{corollary}
\label{cor.ind.gen-com.drop0}Let $S$ be a finite set. Let $a_{s}$ be an
element of $\mathbb{A}$ for each $s\in S$. Let $T$ be a subset of $S$ such
that every $s\in T$ satisfies $a_{s}=0$. Then,%
\[
\sum_{s\in S}a_{s}=\sum_{s\in S\setminus T}a_{s}.
\]

\end{corollary}

\begin{proof}
[Proof of Corollary \ref{cor.ind.gen-com.drop0}.]We have assumed that every
$s\in T$ satisfies $a_{s}=0$. Thus, $\sum_{s\in T}\underbrace{a_{s}}_{=0}%
=\sum_{s\in T}0=0$ (by Theorem \ref{thm.ind.gen-com.sum(0)} (applied to $T$
instead of $S$)).

But $T$ and $S\setminus T$ are subsets of $S$. These two subsets satisfy
$T\cap\left(  S\setminus T\right)  =\varnothing$ and $T\cup\left(  S\setminus
T\right)  =S$ (since $T\subseteq S$). Hence, Theorem
\ref{thm.ind.gen-com.split2} (applied to $X=T$ and $Y=S\setminus T$) yields%
\[
\sum_{s\in S}a_{s}=\underbrace{\sum_{s\in T}a_{s}}_{=0}+\sum_{s\in S\setminus
T}a_{s}=\sum_{s\in S\setminus T}a_{s}.
\]
This proves Corollary \ref{cor.ind.gen-com.drop0}.
\end{proof}

\subsubsection{Substituting the summation index}

Next, we shall show the equality (\ref{eq.sum.subs1}):

\begin{theorem}
\label{thm.ind.gen-com.subst1}Let $S$ and $T$ be two finite sets. Let
$f:S\rightarrow T$ be a \textbf{bijective} map. Let $a_{t}$ be an element of
$\mathbb{A}$ for each $t\in T$. Then,%
\[
\sum_{t\in T}a_{t}=\sum_{s\in S}a_{f\left(  s\right)  }.
\]

\end{theorem}

\begin{proof}
[Proof of Theorem \ref{thm.ind.gen-com.subst1}.]Each $w\in T$ satisfies%
\begin{equation}
\sum_{\substack{s\in S;\\f\left(  s\right)  =w}}a_{f\left(  s\right)  }=a_{w}.
\label{pf.thm.ind.gen-com.subst1.1}%
\end{equation}

[\textit{Proof of (\ref{pf.thm.ind.gen-com.subst1.1}):} Let $w\in T$.

\begin{vershort}
The map $f$ is bijective; thus, it is invertible. In other words, its inverse
map $f^{-1}:T\rightarrow S$ exists. Hence, $f^{-1}\left(  w\right)  $ is a
well-defined element of $S$, and is the only element $t\in S$ satisfying
$f\left(  t\right)  =w$. Therefore,
\begin{equation}
\left\{  t\in S\ \mid\ f\left(  t\right)  =w\right\}  =\left\{  f^{-1}\left(
w\right)  \right\}  . \label{pf.thm.ind.gen-com.subst1.1.pf.short.0}%
\end{equation}

\end{vershort}

\begin{verlong}
The map $f$ is bijective; thus, it is invertible. In other words, its inverse
map $f^{-1}:T\rightarrow S$ exists. Hence, $f^{-1}\left(  w\right)  $ is a
well-defined element of $S$. This element $f^{-1}\left(  w\right)  $ belongs
to $S$ and satisfies $f\left(  f^{-1}\left(  w\right)  \right)  =w$. In other
words, $f^{-1}\left(  w\right)  $ is a $t\in S$ satisfying $f\left(  t\right)
=w$. In other words, $f^{-1}\left(  w\right)  \in\left\{  t\in S\ \mid
\ f\left(  t\right)  =w\right\}  $. Hence,%
\begin{equation}
\left\{  f^{-1}\left(  w\right)  \right\}  \subseteq\left\{  t\in
S\ \mid\ f\left(  t\right)  =w\right\}  .
\label{pf.thm.ind.gen-com.subst1.1.pf.1}%
\end{equation}

On the other hand, let $p\in\left\{  t\in S\ \mid\ f\left(  t\right)
=w\right\}  $. Thus, $p$ is a $t\in S$ satisfying $f\left(  t\right)  =w$. In
other words, $p$ is an element of $S$ and satisfies $f\left(  p\right)  =w$.
From $f\left(  p\right)  =w$, we obtain $p=f^{-1}\left(  w\right)  $ (since
$f^{-1}$ is the inverse of $f$), and thus $p=f^{-1}\left(  w\right)
\in\left\{  f^{-1}\left(  w\right)  \right\}  $. Now, forget that we fixed
$p$. We thus have proven that $p\in\left\{  f^{-1}\left(  w\right)  \right\}
$ for every $p\in\left\{  t\in S\ \mid\ f\left(  t\right)  =w\right\}  $. In
other words,
\[
\left\{  t\in S\ \mid\ f\left(  t\right)  =w\right\}  \subseteq\left\{
f^{-1}\left(  w\right)  \right\}  .
\]
Combining this with (\ref{pf.thm.ind.gen-com.subst1.1.pf.1}), we obtain%
\[
\left\{  t\in S\ \mid\ f\left(  t\right)  =w\right\}  =\left\{  f^{-1}\left(
w\right)  \right\}  .
\]

\end{verlong}

\begin{vershort}
Now,%
\begin{align*}
&  \sum_{\substack{s\in S;\\f\left(  s\right)  =w}}a_{f\left(  s\right)  }\\
&  =\sum_{s\in\left\{  t\in S\ \mid\ f\left(  t\right)  =w\right\}
}a_{f\left(  s\right)  }=\sum_{s\in\left\{  f^{-1}\left(  w\right)  \right\}
}a_{f\left(  s\right)  }\ \ \ \ \ \ \ \ \ \ \left(  \text{by
(\ref{pf.thm.ind.gen-com.subst1.1.pf.short.0})}\right) \\
&  =a_{f\left(  f^{-1}\left(  w\right)  \right)  }\ \ \ \ \ \ \ \ \ \ \left(
\begin{array}
[c]{c}%
\text{by Proposition \ref{prop.ind.gen-com.sum12} \textbf{(a)} (applied to
}\left\{  f^{-1}\left(  w\right)  \right\}  \text{, }a_{f\left(  s\right)  }\\
\text{and }f^{-1}\left(  w\right)  \text{ instead of }S\text{, }a_{s}\text{
and }p\text{)}%
\end{array}
\right) \\
&  =a_{w}\ \ \ \ \ \ \ \ \ \ \left(  \text{since }f\left(  f^{-1}\left(
w\right)  \right)  =w\right)  .
\end{align*}
This proves (\ref{pf.thm.ind.gen-com.subst1.1}).]
\end{vershort}

\begin{verlong}
Now,%
\[
\sum_{\substack{s\in S;\\f\left(  s\right)  =w}}a_{f\left(  s\right)  }%
=\sum_{s\in\left\{  t\in S\ \mid\ f\left(  t\right)  =w\right\}  }a_{f\left(
s\right)  }=\sum_{s\in\left\{  f^{-1}\left(  w\right)  \right\}  }a_{f\left(
s\right)  }%
\]
(since $\left\{  t\in S\ \mid\ f\left(  t\right)  =w\right\}  =\left\{
f^{-1}\left(  w\right)  \right\}  $).

But $\left\{  f^{-1}\left(  w\right)  \right\}  =\left\{  f^{-1}\left(
w\right)  \right\}  $. Hence, Proposition \ref{prop.ind.gen-com.sum12}
\textbf{(a)} (applied to $\left\{  f^{-1}\left(  w\right)  \right\}  $,
$a_{f\left(  s\right)  }$ and $f^{-1}\left(  w\right)  $ instead of $S$,
$a_{s}$ and $p$) yields
\[
\sum_{s\in\left\{  f^{-1}\left(  w\right)  \right\}  }a_{f\left(  s\right)
}=a_{f\left(  f^{-1}\left(  w\right)  \right)  }=a_{w}%
\ \ \ \ \ \ \ \ \ \ \left(  \text{since }f\left(  f^{-1}\left(  w\right)
\right)  =w\right)  .
\]
Hence,%
\[
\sum_{\substack{s\in S;\\f\left(  s\right)  =w}}a_{f\left(  s\right)  }%
=\sum_{s\in\left\{  f^{-1}\left(  w\right)  \right\}  }a_{f\left(  s\right)
}=a_{w}.
\]
This proves (\ref{pf.thm.ind.gen-com.subst1.1}).]
\end{verlong}

Renaming the summation index $w$ as $t$ in the sum $\sum_{w\in T}a_{w}$ does
not change the sum (since $\left(  a_{w}\right)  _{w\in T}$ and $\left(
a_{t}\right)  _{t\in T}$ are the same $\mathbb{A}$-valued $T$-family). In
other words, $\sum_{w\in T}a_{w}=\sum_{t\in T}a_{t}$.

Theorem \ref{thm.ind.gen-com.shephf} (applied to $T$ and $a_{f\left(
s\right)  }$ instead of $W$ and $a_{s}$) yields%
\[
\sum_{s\in S}a_{f\left(  s\right)  }=\sum_{w\in T}\underbrace{\sum
_{\substack{s\in S;\\f\left(  s\right)  =w}}a_{f\left(  s\right)  }%
}_{\substack{=a_{w}\\\text{(by (\ref{pf.thm.ind.gen-com.subst1.1}))}}%
}=\sum_{w\in T}a_{w}=\sum_{t\in T}a_{t}.
\]
This proves Theorem \ref{thm.ind.gen-com.subst1}.
\end{proof}

\subsubsection{Sums of congruences}

Proposition \ref{prop.mod.+-*} \textbf{(a)} says that we can add two
congruences modulo an integer $n$. We shall now see that we can add
\textbf{any} number of congruences modulo an integer $n$:

\begin{theorem}
\label{thm.ind.gen-com.sum-mod}Let $n$ be an integer. Let $S$ be a finite set.
For each $s\in S$, let $a_{s}$ and $b_{s}$ be two integers. Assume that%
\[
a_{s}\equiv b_{s}\operatorname{mod}n\ \ \ \ \ \ \ \ \ \ \text{for each }s\in
S.
\]
Then,%
\[
\sum_{s\in S}a_{s}\equiv\sum_{s\in S}b_{s}\operatorname{mod}n.
\]

\end{theorem}

\begin{proof}
[Proof of Theorem \ref{thm.ind.gen-com.sum-mod}.]We forget that we fixed $n$,
$S$, $a_{s}$ and $b_{s}$. We shall prove Theorem \ref{thm.ind.gen-com.sum-mod}
by induction on $\left\vert S\right\vert $:

\begin{vershort}
\textit{Induction base:} The induction base (i.e., proving that Theorem
\ref{thm.ind.gen-com.sum-mod} holds under the condition that $\left\vert
S\right\vert =0$) is left to the reader (as it boils down to the trivial fact
that $0\equiv0\operatorname{mod}n$).
\end{vershort}

\begin{verlong}
\textit{Induction base:} Theorem \ref{thm.ind.gen-com.sum-mod} holds under the
condition that $\left\vert S\right\vert =0$\ \ \ \ \footnote{\textit{Proof.}
Let $n$, $S$, $a_{s}$ and $b_{s}$ be as in Theorem
\ref{thm.ind.gen-com.sum-mod}. Assume that $\left\vert S\right\vert =0$. Thus,
the first bullet point of Definition \ref{def.ind.gen-com.defsum1} yields
$\sum_{s\in S}a_{s}=0$ and $\sum_{s\in S}b_{s}=0$. Hence,%
\[
\sum_{s\in S}a_{s}=0\equiv0=\sum_{s\in S}b_{s}\operatorname{mod}n.
\]
\par
Now, forget that we fixed $n$, $S$, $a_{s}$ and $b_{s}$. We thus have proved
that if $n$, $S$, $a_{s}$ and $b_{s}$ are as in Theorem
\ref{thm.ind.gen-com.sum-mod}, and if $\left\vert S\right\vert =0$, then
$\sum_{s\in S}a_{s}\equiv\sum_{s\in S}b_{s}\operatorname{mod}n$. In other
words, Theorem \ref{thm.ind.gen-com.sum-mod} holds under the condition that
$\left\vert S\right\vert =0$. Qed.}. This completes the induction base.
\end{verlong}

\textit{Induction step:} Let $m\in\mathbb{N}$. Assume that Theorem
\ref{thm.ind.gen-com.sum-mod} holds under the condition that $\left\vert
S\right\vert =m$. We must now prove that Theorem \ref{thm.ind.gen-com.sum-mod}
holds under the condition that $\left\vert S\right\vert =m+1$.

We have assumed that Theorem \ref{thm.ind.gen-com.sum-mod} holds under the
condition that $\left\vert S\right\vert =m$. In other words, the following
claim holds:

\begin{statement}
\textit{Claim 1:} Let $n$ be an integer. Let $S$ be a finite set such that
$\left\vert S\right\vert =m$. For each $s\in S$, let $a_{s}$ and $b_{s}$ be
two integers. Assume that%
\[
a_{s}\equiv b_{s}\operatorname{mod}n\ \ \ \ \ \ \ \ \ \ \text{for each }s\in
S.
\]
Then,%
\[
\sum_{s\in S}a_{s}\equiv\sum_{s\in S}b_{s}\operatorname{mod}n.
\]

\end{statement}

Next, we shall show the following claim:

\begin{statement}
\textit{Claim 2:} Let $n$ be an integer. Let $S$ be a finite set such that
$\left\vert S\right\vert =m+1$. For each $s\in S$, let $a_{s}$ and $b_{s}$ be
two integers. Assume that%
\begin{equation}
a_{s}\equiv b_{s}\operatorname{mod}n\ \ \ \ \ \ \ \ \ \ \text{for each }s\in
S. \label{pf.thm.ind.gen-com.sum-mod.c2.ass}%
\end{equation}
Then,%
\[
\sum_{s\in S}a_{s}\equiv\sum_{s\in S}b_{s}\operatorname{mod}n.
\]

\end{statement}

[\textit{Proof of Claim 2:} We have $\left\vert S\right\vert =m+1>m\geq0$.
Hence, the set $S$ is nonempty. Thus, there exists some $t\in S$. Consider
this $t$.

From $t\in S$, we obtain $\left\vert S\setminus\left\{  t\right\}  \right\vert
=\left\vert S\right\vert -1=m$ (since $\left\vert S\right\vert =m+1$). Also,
every $s\in S\setminus\left\{  t\right\}  $ satisfies $s\in S\setminus\left\{
t\right\}  \subseteq S$ and thus $a_{s}\equiv b_{s}\operatorname{mod}n$ (by
(\ref{pf.thm.ind.gen-com.sum-mod.c2.ass})). In other words, we have%
\[
a_{s}\equiv b_{s}\operatorname{mod}n\ \ \ \ \ \ \ \ \ \ \text{for each }s\in
S\setminus\left\{  t\right\}  .
\]
Hence, Claim 1 (applied to $S\setminus\left\{  t\right\}  $ instead of $S$)
yields%
\begin{equation}
\sum_{s\in S\setminus\left\{  t\right\}  }a_{s}\equiv\sum_{s\in S\setminus
\left\{  t\right\}  }b_{s}\operatorname{mod}n.
\label{pf.thm.ind.gen-com.sum-mod.c2.pf.1}%
\end{equation}

But $t\in S$. Hence, (\ref{pf.thm.ind.gen-com.sum-mod.c2.ass}) (applied to
$s=t$) yields $a_{t}\equiv b_{t}\operatorname{mod}n$.

Now, Proposition \ref{prop.ind.gen-com.split-off} (applied to $b_{s}$ instead
of $a_{s}$) yields%
\begin{equation}
\sum_{s\in S}b_{s}=b_{t}+\sum_{s\in S\setminus\left\{  t\right\}  }b_{s}.
\label{pf.thm.ind.gen-com.sum-mod.c2.pf.3}%
\end{equation}

But Proposition \ref{prop.ind.gen-com.split-off} yields%
\[
\sum_{s\in S}a_{s}=\underbrace{a_{t}}_{\equiv b_{t}\operatorname{mod}%
n}+\underbrace{\sum_{s\in S\setminus\left\{  t\right\}  }a_{s}}%
_{\substack{\equiv\sum_{s\in S\setminus\left\{  t\right\}  }b_{s}%
\operatorname{mod}n\\\text{(by (\ref{pf.thm.ind.gen-com.sum-mod.c2.pf.1}))}%
}}\equiv b_{t}+\sum_{s\in S\setminus\left\{  t\right\}  }b_{s}=\sum_{s\in
S}b_{s}\operatorname{mod}n
\]
(by (\ref{pf.thm.ind.gen-com.sum-mod.c2.pf.3})). This proves Claim 2.]

But Claim 2 says precisely that Theorem \ref{thm.ind.gen-com.sum-mod} holds
under the condition that $\left\vert S\right\vert =m+1$. Hence, we conclude
that Theorem \ref{thm.ind.gen-com.sum-mod} holds under the condition that
$\left\vert S\right\vert =m+1$ (since Claim 2 is proven). This completes the
induction step. Thus, Theorem \ref{thm.ind.gen-com.sum-mod} is proven by induction.
\end{proof}

As we said, Theorem \ref{thm.ind.gen-com.sum-mod} shows that we can sum up
several congruences. Thus, we can extend our principle of substitutivity for
congruences as follows:

\begin{statement}
\textit{Principle of substitutivity for congruences (stronger version):} Fix
an integer $n$. If two numbers $x$ and $x^{\prime}$ are congruent to each
other modulo $n$ (that is, $x\equiv x^{\prime}\operatorname{mod}n$), and if we
have any expression $A$ that involves only integers, addition, subtraction,
multiplication \textbf{and summation signs}, and involves the object $x$, then
we can replace this $x$ (or, more precisely, any arbitrary appearance of $x$
in $A$) in $A$ by $x^{\prime}$; the value of the resulting expression
$A^{\prime}$ will be congruent to the value of $A$ modulo $n$.
\end{statement}

For example, if $p\in\mathbb{N}$, then%
\[
\sum_{s\in\left\{  1,2,\ldots,p\right\}  }s^{2}\left(  5-3s\right)  \equiv
\sum_{s\in\left\{  1,2,\ldots,p\right\}  }s\left(  5-3s\right)
\operatorname{mod}2
\]
(here, we have replaced the \textquotedblleft$s^{2}$\textquotedblright\ inside
the sum by \textquotedblleft$s$\textquotedblright), because every
$s\in\left\{  1,2,\ldots,p\right\}  $ satisfies $s^{2}\equiv
s\operatorname{mod}2$ (this is easy to check\footnote{\textit{Proof.} Let
$p\in\mathbb{N}$ and $s\in\left\{  1,2,\ldots,p\right\}  $. We must prove that
$s^{2}\equiv s\operatorname{mod}2$.
\par
We have $s\in\left\{  1,2,\ldots,p\right\}  $ and thus $s-1\in\left\{
0,1,\ldots,p-1\right\}  \subseteq\mathbb{N}$. Hence,
(\ref{eq.prop.ind.gen-com.n(n+1)/2.claim}) (applied to $n=s-1$) yields
$\sum_{i\in\left\{  1,2,\ldots,s-1\right\}  }i=\dfrac{\left(  s-1\right)
\left(  \left(  s-1\right)  +1\right)  }{2}=\dfrac{\left(  s-1\right)  s}{2}$.
Hence, $\dfrac{\left(  s-1\right)  s}{2}$ is an integer (since $\sum
_{i\in\left\{  1,2,\ldots,s-1\right\}  }i$ is an integer). In other words,
$2\mid\left(  s-1\right)  s$. In other words, $2\mid s^{2}-s$ (since $\left(
s-1\right)  s=s^{2}-s$). In other words, $s^{2}\equiv s\operatorname{mod}2$
(by the definition of \textquotedblleft congruent\textquotedblright), qed.}).

\subsubsection{\label{subsect.ind.gen-com.prods}Finite products}

Proposition \ref{prop.ind.gen-com.fgh} is a property of the addition of
numbers; it has an analogue for multiplication of numbers:

\begin{proposition}
\label{prop.ind.gen-com.fgh*}Let $a$, $b$ and $c$ be three numbers (i.e.,
elements of $\mathbb{A}$). Then, $\left(  ab\right)  c=a\left(  bc\right)  $.
\end{proposition}

Proposition \ref{prop.ind.gen-com.fgh*} is known as the \textit{associativity
of multiplication} (in $\mathbb{A}$), and is fundamental; its proof can be
found in any textbook on the construction of the number system\footnote{For
example, Proposition \ref{prop.ind.gen-com.fgh*} is proven in \cite[Theorem
3.2.3 (7)]{Swanso18} for the case when $\mathbb{A}=\mathbb{N}$; in
\cite[Theorem 3.5.4 (7)]{Swanso18} for the case when $\mathbb{A}=\mathbb{Z}$;
in \cite[Theorem 3.6.4 (7)]{Swanso18} for the case when $\mathbb{A}%
=\mathbb{Q}$; in \cite[Theorem 3.7.14]{Swanso18} for the case when
$\mathbb{A}=\mathbb{R}$; in \cite[Theorem 3.9.2]{Swanso18} for the case when
$\mathbb{A}=\mathbb{C}$.}.

Proposition \ref{prop.ind.gen-com.fg} also has an analogue for multiplication:

\begin{proposition}
\label{prop.ind.gen-com.fg*}Let $a$ and $b$ be two numbers (i.e., elements of
$\mathbb{A}$). Then, $ab=ba$.
\end{proposition}

Proposition \ref{prop.ind.gen-com.fg*} is known as the \textit{commutativity
of multiplication} (in $\mathbb{A}$), and again is a fundamental result whose
proofs are found in standard textbooks\footnote{For example, Proposition
\ref{prop.ind.gen-com.fg*} is proven in \cite[Theorem 3.2.3 (8)]{Swanso18} for
the case when $\mathbb{A}=\mathbb{N}$; in \cite[Theorem 3.5.4 (8)]{Swanso18}
for the case when $\mathbb{A}=\mathbb{Z}$; in \cite[Theorem 3.6.4
(8)]{Swanso18} for the case when $\mathbb{A}=\mathbb{Q}$; in \cite[Theorem
3.7.14]{Swanso18} for the case when $\mathbb{A}=\mathbb{R}$; in \cite[Theorem
3.9.2]{Swanso18} for the case when $\mathbb{A}=\mathbb{C}$.}.

Proposition \ref{prop.ind.gen-com.distr} has an analogue for multiplication as
well (but note that $x$ now needs to be in $\mathbb{N}$, in order to guarantee
that the powers are well-defined):

\begin{proposition}
\label{prop.ind.gen-com.distr*}Let $x\in\mathbb{N}$. Let $y$ and $z$ be two
numbers (i.e., elements of $\mathbb{A}$). Then, $\left(  yz\right)  ^{x}%
=y^{x}z^{x}$.
\end{proposition}

Proposition \ref{prop.ind.gen-com.distr*} is one of the laws of exponents, and
can easily be shown by induction on $x$ (using Proposition
\ref{prop.ind.gen-com.fg*} and Proposition \ref{prop.ind.gen-com.fgh*}).

So far in Section \ref{sect.ind.gen-com}, we have been studying \textbf{sums}
of $\mathbb{A}$-valued $S$-families (when $S$ is a finite set): We have proven
that the definition of $\sum_{s\in S}a_{s}$ given in Section
\ref{sect.sums-repetitorium} is legitimate, and we have proven several
properties of such sums. By the exact same reasoning (but with addition
replaced by multiplication), we can study \textbf{products} of $\mathbb{A}%
$-valued $S$-families. In particular, we can similarly prove that the
definition of $\prod_{s\in S}a_{s}$ given in Section
\ref{sect.sums-repetitorium} is legitimate, and we can prove properties of
such products that are analogous to the properties of sums proven above
(except for Proposition \ref{prop.ind.gen-com.n(n+1)/2}, which does not have
an analogue for products)\footnote{We need to be slightly careful when we
adapt our above proofs to products instead of sums: Apart from replacing
addition by multiplication everywhere, we need to:
\par
\begin{itemize}
\item replace the number $0$ by $1$ whenever it appears in a computation
inside $\mathbb{A}$ (but, of course, not when it appears as the size of a
set);
\par
\item replace every $\sum$ sign by a $\prod$ sign;
\par
\item replace \textquotedblleft let $\lambda$ be an element of $\mathbb{A}%
$\textquotedblright\ by \textquotedblleft let $\lambda$ be an element of
$\mathbb{N}$\textquotedblright\ in Theorem \ref{thm.ind.gen-com.sum(la)};
\par
\item replace any expression of the form \textquotedblleft$\lambda
b$\textquotedblright\ by \textquotedblleft$b^{\lambda}$\textquotedblright\ in
Theorem \ref{thm.ind.gen-com.sum(la)} (so that the claim of Theorem
\ref{thm.ind.gen-com.sum(la)} becomes $\prod_{s\in S}\left(  a_{s}\right)
^{\lambda}=\left(  \prod_{s\in S}a_{s}\right)  ^{\lambda}$) and in its proof;
\par
\item replace every reference to Proposition \ref{prop.ind.gen-com.fgh} by a
reference to Proposition \ref{prop.ind.gen-com.fgh*};
\par
\item replace every reference to Proposition \ref{prop.ind.gen-com.fg} by a
reference to Proposition \ref{prop.ind.gen-com.fg*};
\par
\item replace every reference to Proposition \ref{prop.ind.gen-com.distr} by a
reference to Proposition \ref{prop.ind.gen-com.distr*}.
\end{itemize}
\par
And, to be fully precise: We should not replace addition by multiplication
\textbf{everywhere} (e.g., we should not replace \textquotedblleft$\left\vert
S\right\vert =m+1$\textquotedblright\ by \textquotedblleft$\left\vert
S\right\vert =m\cdot1$\textquotedblright\ in the proof of Theorem
\ref{thm.ind.gen-com.shephf}), but of course only where it stands for the
addition \textbf{inside }$\mathbb{A}$.}. For example, the following theorems
are analogues of Theorem \ref{thm.ind.gen-com.sum(a+b)}, Theorem
\ref{thm.ind.gen-com.sum(la)}, Theorem \ref{thm.ind.gen-com.sum(0)}, Theorem
\ref{thm.ind.gen-com.shephf}, Theorem \ref{thm.ind.gen-com.subst1} and Theorem
\ref{thm.ind.gen-com.sum-mod}, respectively:

\begin{theorem}
\label{thm.ind.gen-com.prod(a+b)}Let $S$ be a finite set. For every $s\in S$,
let $a_{s}$ and $b_{s}$ be elements of $\mathbb{A}$. Then,%
\[
\prod_{s\in S}\left(  a_{s}b_{s}\right)  =\left(  \prod_{s\in S}a_{s}\right)
\cdot\left(  \prod_{s\in S}b_{s}\right)  .
\]

\end{theorem}

\begin{theorem}
\label{thm.ind.gen-com.prod(la)}Let $S$ be a finite set. For every $s\in S$,
let $a_{s}$ be an element of $\mathbb{A}$. Also, let $\lambda$ be an element
of $\mathbb{N}$. Then,%
\[
\prod_{s\in S}\left(  a_{s}\right)  ^{\lambda}=\left(  \prod_{s\in S}%
a_{s}\right)  ^{\lambda}.
\]

\end{theorem}

\begin{theorem}
\label{thm.ind.gen-com.prod(1)}Let $S$ be a finite set. Then,
\[
\prod_{s\in S}1=1.
\]

\end{theorem}

\begin{theorem}
\label{thm.ind.gen-com.shephf*}Let $S$ be a finite set. Let $W$ be a finite
set. Let $f:S\rightarrow W$ be a map. Let $a_{s}$ be an element of
$\mathbb{A}$ for each $s\in S$. Then,%
\[
\prod_{s\in S}a_{s}=\prod_{w\in W}\prod_{\substack{s\in S;\\f\left(  s\right)
=w}}a_{s}.
\]

\end{theorem}

\begin{theorem}
\label{thm.ind.gen-com.subst1*}Let $S$ and $T$ be two finite sets. Let
$f:S\rightarrow T$ be a \textbf{bijective} map. Let $a_{t}$ be an element of
$\mathbb{A}$ for each $t\in T$. Then,%
\[
\prod_{t\in T}a_{t}=\prod_{s\in S}a_{f\left(  s\right)  }.
\]

\end{theorem}

\begin{theorem}
\label{thm.ind.gen-com.prod-mod}Let $n$ be an integer. Let $S$ be a finite
set. For each $s\in S$, let $a_{s}$ and $b_{s}$ be two integers. Assume that%
\[
a_{s}\equiv b_{s}\operatorname{mod}n\ \ \ \ \ \ \ \ \ \ \text{for each }s\in
S.
\]
Then,%
\[
\prod_{s\in S}a_{s}\equiv\prod_{s\in S}b_{s}\operatorname{mod}n.
\]

\end{theorem}

\subsubsection{Finitely supported (but possibly infinite) sums}

In Section \ref{sect.sums-repetitorium}, we mentioned that a sum of the form
$\sum_{s\in S}a_{s}$ can be well-defined even when the set $S$ is not finite.
Indeed, for it to be well-defined, it suffices that \textbf{only finitely many
among the }$a_{s}$ \textbf{are nonzero} (or, more rigorously: only finitely
many $s\in S$ satisfy $a_{s}\neq0$). As we already mentioned, the sum
$\sum_{s\in S}a_{s}$ in this case is defined by discarding the zero addends
and summing the finitely many addends that remain. Let us briefly discuss such
sums (without focussing on advanced properties):

\begin{definition}
\label{def.ind.gen-com.fin-sup.def}Let $S$ be any set. An $\mathbb{A}$-valued
$S$-family $\left(  a_{s}\right)  _{s\in S}$ is said to be \textit{finitely
supported} if only finitely many $s\in S$ satisfy $a_{s}\neq0$.
\end{definition}

So the sums we want to discuss are sums $\sum_{s\in S}a_{s}$ for which the set
$S$ may be infinite but the $S$-family $\left(  a_{s}\right)  _{s\in S}$ is
finitely supported. Let us repeat the definition of such sums in more rigorous language:

\begin{definition}
\label{def.ind.gen-com.fin-sup.sum}Let $S$ be any set. Let $\left(
a_{s}\right)  _{s\in S}$ be a finitely supported $\mathbb{A}$-valued
$S$-family. Thus, there exists a \textbf{finite} subset $T$ of $S$ such that%
\begin{equation}
\text{every }s\in S\setminus T\text{ satisfies }a_{s}=0.
\label{eq.def.ind.gen-com.fin-sup.sum.0}%
\end{equation}
(This is because only finitely many $s\in S$ satisfy $a_{s}\neq0$.) We then
define the sum $\sum_{s\in S}a_{s}$ to be $\sum_{s\in T}a_{s}$. (This
definition is legitimate, because Proposition
\ref{prop.ind.gen-com.fin-sup.leg} \textbf{(a)} below shows that $\sum_{s\in
T}a_{s}$ does not depend on the choice of $T$.)
\end{definition}

This definition formalizes what we said above about making sense of
$\sum_{s\in S}a_{s}$: Namely, we discard zero addends (namely, the addends
corresponding to $s\in S\setminus T$) and only sum the finitely many addends
that remain (these are the addends corresponding to $s\in T$); thus, we get
$\sum_{s\in T}a_{s}$. Note that we are not requiring that every $s\in T$
satisfies $a_{s}\neq0$; that is, we are not necessarily discarding
\textbf{all} the zero addends from our sum (but merely discarding enough of
them to ensure that only finitely many remain). This may appear like a strange
choice (why introduce extra freedom into the definition?), but is reasonable
from the viewpoint of constructive mathematics (where it is not always
decidable if a number is $0$ or not).

\begin{proposition}
\label{prop.ind.gen-com.fin-sup.leg}Let $S$ be any set. Let $\left(
a_{s}\right)  _{s\in S}$ be a finitely supported $\mathbb{A}$-valued $S$-family.

\textbf{(a)} If $T$ is a finite subset of $S$ such that
(\ref{eq.def.ind.gen-com.fin-sup.sum.0}) holds, then the sum $\sum_{s\in
T}a_{s}$ does not depend on the choice of $T$. (That is, if $T_{1}$ and
$T_{2}$ are two finite subsets $T$ of $S$ satisfying
(\ref{eq.def.ind.gen-com.fin-sup.sum.0}), then $\sum_{s\in T_{1}}a_{s}%
=\sum_{s\in T_{2}}a_{s}$.)

\textbf{(b)} If the set $S$ is finite, then the sum $\sum_{s\in S}a_{s}$
defined in Definition \ref{def.ind.gen-com.fin-sup.sum} is identical with the
sum $\sum_{s\in S}a_{s}$ defined in Definition \ref{def.ind.gen-com.defsum1}.
(Thus, Definition \ref{def.ind.gen-com.fin-sup.sum} does not conflict with the
previous definition of $\sum_{s\in S}a_{s}$ for finite sets $S$.)
\end{proposition}

Proposition \ref{prop.ind.gen-com.fin-sup.leg} is fairly easy to prove using
Corollary \ref{cor.ind.gen-com.drop0}; this proof is part of Exercise
\ref{exe.ind.gen-com.fin-sup.proofs} below.

Most properties of finite sums have analogues for sums of finitely supported
$\mathbb{A}$-valued $S$-families. For example, here is an analogue of Theorem
\ref{thm.ind.gen-com.sum(a+b)}:

\begin{theorem}
\label{thm.ind.gen-com.sum(a+b).gen}Let $S$ be a set. Let $\left(
a_{s}\right)  _{s\in S}$ and $\left(  b_{s}\right)  _{s\in S}$ be two finitely
supported $\mathbb{A}$-valued $S$-families. Then, the $\mathbb{A}$-valued
$S$-family $\left(  a_{s}+b_{s}\right)  _{s\in S}$ is finitely supported as
well, and we have%
\[
\sum_{s\in S}\left(  a_{s}+b_{s}\right)  =\sum_{s\in S}a_{s}+\sum_{s\in
S}b_{s}.
\]

\end{theorem}

The proof of Theorem \ref{thm.ind.gen-com.sum(a+b).gen} is fairly simple (it
relies prominently on the fact that the union of two finite sets is finite),
and again is part of Exercise \ref{exe.ind.gen-com.fin-sup.proofs} below.

It is also easy to state and prove analogues of Theorem
\ref{thm.ind.gen-com.sum(la)} and Theorem \ref{thm.ind.gen-com.sum(0)}.

We can next prove (\ref{eq.sum.sheph}) in full generality (not only when $W$
is finite):

\begin{theorem}
\label{thm.ind.gen-com.sheph}Let $S$ be a finite set. Let $W$ be a set. Let
$f:S\rightarrow W$ be a map. Let $a_{s}$ be an element of $\mathbb{A}$ for
each $s\in S$. Then, the $\mathbb{A}$-valued $W$-family $\left(
\sum_{\substack{s\in S;\\f\left(  s\right)  =w}}a_{s}\right)  _{w\in W}$ is
finitely supported and satisfies%
\[
\sum_{s\in S}a_{s}=\sum_{w\in W}\sum_{\substack{s\in S;\\f\left(  s\right)
=w}}a_{s}.
\]

\end{theorem}

Note that the sum on the right hand side of Theorem
\ref{thm.ind.gen-com.sheph} makes sense even when $W$ is infinite, because the
$W$-family $\left(  \sum_{\substack{s\in S;\\f\left(  s\right)  =w}%
}a_{s}\right)  _{w\in W}$ is finitely supported (i.e., only finitely many
$w\in W$ satisfy $\sum_{\substack{s\in S;\\f\left(  s\right)  =w}}a_{s}\neq
0$). The easiest way to prove Theorem \ref{thm.ind.gen-com.sheph} is probably
by reducing it to Theorem \ref{thm.ind.gen-com.shephf} (since $f\left(
S\right)  $ is a finite subset of $W$, and every $w\in W\setminus f\left(
S\right)  $ satisfies $\sum_{\substack{s\in S;\\f\left(  s\right)  =w}%
}a_{s}=\left(  \text{empty sum}\right)  =0$). Again, we leave the details to
the interested reader.

Again, we refer to Exercise \ref{exe.ind.gen-com.fin-sup.proofs} for the proof
of Theorem \ref{thm.ind.gen-com.sheph}.

Actually, Theorem \ref{thm.ind.gen-com.sheph} can be generalized even further:

\begin{theorem}
\label{thm.ind.gen-com.sheph-gen}Let $S$ be a set. Let $W$ be a set. Let
$f:S\rightarrow W$ be a map. Let $\left(  a_{s}\right)  _{s\in S}$ be a
finitely supported $\mathbb{A}$-valued $S$-family. Then, for each $w\in W$,
the $\mathbb{A}$-valued $\left\{  t\in S\ \mid\ f\left(  t\right)  =w\right\}
$-family $\left(  a_{s}\right)  _{s\in\left\{  t\in S\ \mid\ f\left(
t\right)  =w\right\}  }$ is finitely supported as well (so that the sum
$\sum_{\substack{s\in S;\\f\left(  s\right)  =w}}a_{s}$ is well-defined).
Furthermore, the $\mathbb{A}$-valued $W$-family $\left(  \sum_{\substack{s\in
S;\\f\left(  s\right)  =w}}a_{s}\right)  _{w\in W}$ is also finitely
supported. Finally,%
\[
\sum_{s\in S}a_{s}=\sum_{w\in W}\sum_{\substack{s\in S;\\f\left(  s\right)
=w}}a_{s}.
\]

\end{theorem}

Again, see Exercise \ref{exe.ind.gen-com.fin-sup.proofs} for the proof. This
theorem can be used to obtain an analogue of Theorem
\ref{thm.ind.gen-com.split2} for finitely supported $\mathbb{A}$-valued $S$-families.

\begin{exercise}
\label{exe.ind.gen-com.fin-sup.proofs}Prove Proposition
\ref{prop.ind.gen-com.fin-sup.leg}, Theorem \ref{thm.ind.gen-com.sum(a+b).gen}%
, Theorem \ref{thm.ind.gen-com.sheph} and Theorem
\ref{thm.ind.gen-com.sheph-gen}.
\end{exercise}

Thus, we have defined the values of certain infinite sums (although not nearly
as many infinite sums as analysis can make sense of). We can similarly define
the values of certain infinite products: In order for $\prod_{s\in S}a_{s}$ to
be well-defined, it suffices that \textbf{only finitely many among the }%
$a_{s}$ \textbf{are distinct from }$1$ (or, more rigorously: only finitely
many $s\in S$ satisfy $a_{s}\neq1$). We leave the details and properties of
this definition to the reader.

\subsection{Two-sided induction}

\subsubsection{The principle of two-sided induction}

Let us now return to studying induction principles. We have seen several
induction principles that allow us to prove statements about nonnegative
integers, integers in $\mathbb{Z}_{\geq g}$ or integers in an interval. What
about proving statements about \textbf{arbitrary} integers? The induction
principles we have seen so far do not suffice to prove such statements
directly, since our induction steps always \textquotedblleft go
up\textquotedblright\ (in the sense that they begin by assuming that our
statement $\mathcal{A}\left(  k\right)  $ holds for some integers $k$, and
involve proving that it also holds for a \textbf{larger} value of $k$), but it
is impossible to traverse all the integers by starting at some integer $g$ and
going up (you will never get to $g-1$ this way). In contrast, the following
induction principle includes both an \textquotedblleft
upwards\textquotedblright\ and a \textquotedblleft downwards\textquotedblright%
\ induction step, which makes it suited for proving statements about all integers:

\begin{theorem}
\label{thm.ind.IPg+-}Let $g\in\mathbb{Z}$. Let $\mathbb{Z}_{\leq g}$ be the
set $\left\{  g,g-1,g-2,\ldots\right\}  $ (that is, the set of all integers
that are $\leq g$).

For each $n\in\mathbb{Z}$, let $\mathcal{A}\left(  n\right)  $ be a logical statement.

Assume the following:

\begin{statement}
\textit{Assumption 1:} The statement $\mathcal{A}\left(  g\right)  $ holds.
\end{statement}

\begin{statement}
\textit{Assumption 2:} If $m\in\mathbb{Z}_{\geq g}$ is such that
$\mathcal{A}\left(  m\right)  $ holds, then $\mathcal{A}\left(  m+1\right)  $
also holds.
\end{statement}

\begin{statement}
\textit{Assumption 3:} If $m\in\mathbb{Z}_{\leq g}$ is such that
$\mathcal{A}\left(  m\right)  $ holds, then $\mathcal{A}\left(  m-1\right)  $
also holds.
\end{statement}

Then, $\mathcal{A}\left(  n\right)  $ holds for each $n\in\mathbb{Z}$.
\end{theorem}

Theorem \ref{thm.ind.IPg+-} is known as the \textit{principle of two-sided
induction}. Roughly speaking, a proof using Theorem \ref{thm.ind.IPg+-} will
involve two induction steps: one that \textquotedblleft goes
up\textquotedblright\ (proving that Assumption 2 holds), and one that
\textquotedblleft goes down\textquotedblright\ (proving that Assumption 3
holds). However, in practice, Theorem \ref{thm.ind.IPg+-} is seldom used,
which is why we shall not make any conventions about how to write proofs using
Theorem \ref{thm.ind.IPg+-}. We will only give one example for such a proof.

Let us first prove Theorem \ref{thm.ind.IPg+-} itself:

\begin{proof}
[Proof of Theorem \ref{thm.ind.IPg+-}.]Assumptions 1 and 2 of Theorem
\ref{thm.ind.IPg+-} are exactly Assumptions 1 and 2 of Theorem
\ref{thm.ind.IPg}. Hence, Assumptions 1 and 2 of Theorem \ref{thm.ind.IPg}
hold (since Assumptions 1 and 2 of Theorem \ref{thm.ind.IPg+-} hold). Thus,
Theorem \ref{thm.ind.IPg} shows that%
\begin{equation}
\mathcal{A}\left(  n\right)  \text{ holds for each }n\in\mathbb{Z}_{\geq g}.
\label{pf.thm.ind.IPg+-.pospart}%
\end{equation}

On the other hand, for each $n\in\mathbb{Z}$, we define a logical statement
$\mathcal{B}\left(  n\right)  $ by $\mathcal{B}\left(  n\right)
=\mathcal{A}\left(  2g-n\right)  $. We shall now consider the Assumptions A
and B of Corollary \ref{cor.ind.IPg.renamed}.

The definition of $\mathcal{B}\left(  g\right)  $ yields $\mathcal{B}\left(
g\right)  =\mathcal{A}\left(  2g-g\right)  =\mathcal{A}\left(  g\right)  $
(since $2g-g=g$). Hence, the statement $\mathcal{B}\left(  g\right)  $ holds
(since the statement $\mathcal{A}\left(  g\right)  $ holds (by Assumption 1)).
In other words, Assumption A is satisfied.

Next, let $p\in\mathbb{Z}_{\geq g}$ be such that $\mathcal{B}\left(  p\right)
$ holds. We shall show that $\mathcal{B}\left(  p+1\right)  $ holds.

Indeed, we have $\mathcal{B}\left(  p\right)  =\mathcal{A}\left(  2g-p\right)
$ (by the definition of $\mathcal{B}\left(  p\right)  $). Thus, $\mathcal{A}%
\left(  2g-p\right)  $ holds (since $\mathcal{B}\left(  p\right)  $ holds).
But $p\in\mathbb{Z}_{\geq g}$; hence, $p$ is an integer that is $\geq g$.
Thus, $p\geq g$, so that $2g-\underbrace{p}_{\geq g}\leq2g-g=g$. Hence, $2g-p$
is an integer that is $\leq g$. In other words, $2g-p\in\mathbb{Z}_{\leq g}$.
Therefore, Assumption 3 (applied to $m=2g-p$) shows that $\mathcal{A}\left(
2g-p-1\right)  $ also holds (since $\mathcal{A}\left(  2g-p\right)  $ holds).
But the definition of $\mathcal{B}\left(  p+1\right)  $ yields $\mathcal{B}%
\left(  p+1\right)  =\mathcal{A}\left(  \underbrace{2g-\left(  p+1\right)
}_{=2g-p-1}\right)  =\mathcal{A}\left(  2g-p-1\right)  $. Hence,
$\mathcal{B}\left(  p+1\right)  $ holds (since $\mathcal{A}\left(
2g-p-1\right)  $ holds).

Now, forget that we fixed $p$. We thus have shown that if $p\in\mathbb{Z}%
_{\geq g}$ is such that $\mathcal{B}\left(  p\right)  $ holds, then
$\mathcal{B}\left(  p+1\right)  $ also holds. In other words, Assumption B is satisfied.

We now have shown that both Assumptions A and B are satisfied. Hence,
Corollary \ref{cor.ind.IPg.renamed} shows that%
\begin{equation}
\mathcal{B}\left(  n\right)  \text{ holds for each }n\in\mathbb{Z}_{\geq g}.
\label{pf.thm.ind.IPg+-.negpart1}%
\end{equation}

Now, let $n\in\mathbb{Z}$. We shall prove that $\mathcal{A}\left(  n\right)  $ holds.

Indeed, we have either $n\geq g$ or $n<g$. Hence, we are in one of the
following two cases:

\textit{Case 1:} We have $n\geq g$.

\textit{Case 2:} We have $n<g$.

Let us first consider Case 1. In this case, we have $n\geq g$. Hence,
$n\in\mathbb{Z}_{\geq g}$ (since $n$ is an integer). Thus,
(\ref{pf.thm.ind.IPg+-.pospart}) shows that $\mathcal{A}\left(  n\right)  $
holds. We thus have proven that $\mathcal{A}\left(  n\right)  $ holds in Case 1.

Let us now consider Case 2. In this case, we have $n<g$. Thus, $n\leq g$.
Hence, $2g-\underbrace{n}_{\leq g}\geq2g-g=g$. Thus, $2g-n\in\mathbb{Z}_{\geq
g}$ (since $2g-n$ is an integer). Hence, (\ref{pf.thm.ind.IPg+-.negpart1})
(applied to $2g-n$ instead of $n$) shows that $\mathcal{B}\left(  2g-n\right)
$ holds. But the definition of $\mathcal{B}\left(  2g-n\right)  $ yields
$\mathcal{B}\left(  2g-n\right)  =\mathcal{A}\left(  \underbrace{2g-\left(
2g-n\right)  }_{=n}\right)  =\mathcal{A}\left(  n\right)  $. Hence,
$\mathcal{A}\left(  n\right)  $ holds (since $\mathcal{B}\left(  2g-n\right)
$ holds). Thus, we have proven that $\mathcal{A}\left(  n\right)  $ holds in
Case 2.

We now have shown that $\mathcal{A}\left(  n\right)  $ holds in each of the
two Cases 1 and 2. Since these two Cases cover all possibilities, we thus
conclude that $\mathcal{A}\left(  n\right)  $ always holds.

Now, forget that we fixed $n$. We thus have proven that $\mathcal{A}\left(
n\right)  $ holds for each $n\in\mathbb{Z}$. This proves Theorem
\ref{thm.ind.IPg+-}.
\end{proof}

As an example for the use of Theorem \ref{thm.ind.IPg+-}, we shall prove the
following fact:

\begin{proposition}
\label{prop.ind.quo-rem-ex}Let $N$ be a positive integer. For each
$n\in\mathbb{Z}$, there exist $q\in\mathbb{Z}$ and $r\in\left\{
0,1,\ldots,N-1\right\}  $ such that $n=qN+r$.
\end{proposition}

We shall soon see (in Theorem \ref{thm.ind.quo-rem}) that these $q$ and $r$
are actually uniquely determined by $N$ and $n$; they are called the
\textit{quotient} and the \textit{remainder of the division of }$n$ \textit{by
}$N$. This is fundamental to all of number theory.

\begin{example}
If we apply Proposition \ref{prop.ind.quo-rem-ex} to $N=4$ and $n=10$, then we
conclude that there exist $q\in\mathbb{Z}$ and $r\in\left\{  0,1,2,3\right\}
$ such that $10=q\cdot4+r$. And indeed, such $q$ and $r$ can easily be found
($q=2$ and $r=2$).
\end{example}

\begin{proof}
[Proof of Proposition \ref{prop.ind.quo-rem-ex}.]First, we notice that
$N-1\in\mathbb{N}$ (since $N$ is a positive integer). Hence, $0\in\left\{
0,1,\ldots,N-1\right\}  $ and $N-1\in\left\{  0,1,\ldots,N-1\right\}  $.

For each $n\in\mathbb{Z}$, we let $\mathcal{A}\left(  n\right)  $ be the
statement%
\[
\left(  \text{there exist }q\in\mathbb{Z}\text{ and }r\in\left\{
0,1,\ldots,N-1\right\}  \text{ such that }n=qN+r\right)  .
\]

Let $g=0$; thus, $g\in\mathbb{Z}$. Define the set $\mathbb{Z}_{\leq g}$ as in
Theorem \ref{thm.ind.IPg+-}. (Thus, $\mathbb{Z}_{\leq g}=\left\{
g,g-1,g-2,\ldots\right\}  =\left\{  0,-1,-2,\ldots\right\}  $ is the set of
all nonpositive integers.) We shall now show that Assumptions 1, 2 and 3 of
Theorem \ref{thm.ind.IPg+-} are satisfied.

[\textit{Proof that Assumption 1 is satisfied:} We have $0\in\left\{
0,1,\ldots,N-1\right\}  $ and $0=0N+0$. Hence, there exist $q\in\mathbb{Z}$
and $r\in\left\{  0,1,\ldots,N-1\right\}  $ such that $0=qN+r$ (namely, $q=0$
and $r=0$). In other words, the statement $\mathcal{A}\left(  0\right)  $
holds\footnote{since the statement $\mathcal{A}\left(  0\right)  $ is defined
as \newline$\left(  \text{there exist }q\in\mathbb{Z}\text{ and }r\in\left\{
0,1,\ldots,N-1\right\}  \text{ such that }0=qN+r\right)  $}. In other words,
the statement $\mathcal{A}\left(  g\right)  $ holds (since $g=0$). In other
words, Assumption 1 is satisfied.]

[\textit{Proof that Assumption 2 is satisfied:} Let $m\in\mathbb{Z}_{\geq g}$
be such that $\mathcal{A}\left(  m\right)  $ holds. We shall show that
$\mathcal{A}\left(  m+1\right)  $ also holds.

We know that $\mathcal{A}\left(  m\right)  $ holds. In other words, there
exist $q\in\mathbb{Z}$ and $r\in\left\{  0,1,\ldots,N-1\right\}  $ such that
$m=qN+r$\ \ \ \ \footnote{since the statement $\mathcal{A}\left(  m\right)  $
is defined as \newline$\left(  \text{there exist }q\in\mathbb{Z}\text{ and
}r\in\left\{  0,1,\ldots,N-1\right\}  \text{ such that }m=qN+r\right)  $}.
Consider these $q$ and $r$, and denote them by $q_{0}$ and $r_{0}$. Thus,
$q_{0}\in\mathbb{Z}$ and $r_{0}\in\left\{  0,1,\ldots,N-1\right\}  $ and
$m=q_{0}N+r_{0}$.

Now, we are in one of the following two cases:

\textit{Case 1:} We have $r_{0}=N-1$.

\textit{Case 2:} We have $r_{0}\neq N-1$.

Let us first consider Case 1. In this case, we have $r_{0}=N-1$. Hence,
$r_{0}+1=N$. Now,%
\[
\underbrace{m}_{=q_{0}N+r_{0}}+1=q_{0}N+\underbrace{r_{0}+1}_{=N}%
=q_{0}N+N=\left(  q_{0}+1\right)  N=\left(  q_{0}+1\right)  N+0.
\]
Since $0\in\left\{  0,1,\ldots,N-1\right\}  $, we can therefore conclude that
there exist $q\in\mathbb{Z}$ and $r\in\left\{  0,1,\ldots,N-1\right\}  $ such
that $m+1=qN+r$ (namely, $q=q_{0}+1$ and $r=0$). In other words, the statement
$\mathcal{A}\left(  m+1\right)  $ holds\footnote{since the statement
$\mathcal{A}\left(  m+1\right)  $ is defined as \newline$\left(  \text{there
exist }q\in\mathbb{Z}\text{ and }r\in\left\{  0,1,\ldots,N-1\right\}  \text{
such that }m+1=qN+r\right)  $}. Thus, we have proven that $\mathcal{A}\left(
m+1\right)  $ holds in Case 1.

Let us next consider Case 2. In this case, we have $r_{0}\neq N-1$. Combining
$r_{0}\in\left\{  0,1,\ldots,N-1\right\}  $ with $r_{0}\neq N-1$, we obtain%
\[
r_{0}\in\left\{  0,1,\ldots,N-1\right\}  \setminus\left\{  N-1\right\}
=\left\{  0,1,\ldots,N-2\right\}  ,
\]
so that $r_{0}+1\in\left\{  1,2,\ldots,N-1\right\}  \subseteq\left\{
0,1,\ldots,N-1\right\}  $. Also,%
\[
\underbrace{m}_{=q_{0}N+r_{0}}+1=q_{0}N+r_{0}+1=q_{0}N+\left(  r_{0}+1\right)
.
\]
Since $r_{0}+1\in\left\{  0,1,\ldots,N-1\right\}  $, we can therefore conclude
that there exist $q\in\mathbb{Z}$ and $r\in\left\{  0,1,\ldots,N-1\right\}  $
such that $m+1=qN+r$ (namely, $q=q_{0}$ and $r=r_{0}+1$). In other words, the
statement $\mathcal{A}\left(  m+1\right)  $ holds\footnote{since the statement
$\mathcal{A}\left(  m+1\right)  $ is defined as \newline$\left(  \text{there
exist }q\in\mathbb{Z}\text{ and }r\in\left\{  0,1,\ldots,N-1\right\}  \text{
such that }m+1=qN+r\right)  $}. Thus, we have proven that $\mathcal{A}\left(
m+1\right)  $ holds in Case 2.

We have now proven that $\mathcal{A}\left(  m+1\right)  $ holds in each of the
two Cases 1 and 2. Since these two Cases cover all possibilities, we thus
conclude that $\mathcal{A}\left(  m+1\right)  $ always holds.

Now, forget that we fixed $m$. We thus have shown that if $m\in\mathbb{Z}%
_{\geq g}$ is such that $\mathcal{A}\left(  m\right)  $ holds, then
$\mathcal{A}\left(  m+1\right)  $ also holds. In other words, Assumption 2 is satisfied.]

[\textit{Proof that Assumption 3 is satisfied:} Let $m\in\mathbb{Z}_{\leq g}$
be such that $\mathcal{A}\left(  m\right)  $ holds. We shall show that
$\mathcal{A}\left(  m-1\right)  $ also holds.

We know that $\mathcal{A}\left(  m\right)  $ holds. In other words, there
exist $q\in\mathbb{Z}$ and $r\in\left\{  0,1,\ldots,N-1\right\}  $ such that
$m=qN+r$\ \ \ \ \footnote{since the statement $\mathcal{A}\left(  m\right)  $
is defined as \newline$\left(  \text{there exist }q\in\mathbb{Z}\text{ and
}r\in\left\{  0,1,\ldots,N-1\right\}  \text{ such that }m=qN+r\right)  $}.
Consider these $q$ and $r$, and denote them by $q_{0}$ and $r_{0}$. Thus,
$q_{0}\in\mathbb{Z}$ and $r_{0}\in\left\{  0,1,\ldots,N-1\right\}  $ and
$m=q_{0}N+r_{0}$.

Now, we are in one of the following two cases:

\textit{Case 1:} We have $r_{0}=0$.

\textit{Case 2:} We have $r_{0}\neq0$.

Let us first consider Case 1. In this case, we have $r_{0}=0$. Now,%
\begin{align*}
\underbrace{m}_{=q_{0}N+r_{0}}-1  &  =q_{0}N+\underbrace{r_{0}}_{=0}%
-1=q_{0}N-1=\underbrace{q_{0}N-N}_{=\left(  q_{0}-1\right)  N}+\left(
N-1\right) \\
&  =\left(  q_{0}-1\right)  N+\left(  N-1\right)  .
\end{align*}
Since $N-1\in\left\{  0,1,\ldots,N-1\right\}  $, we can therefore conclude
that there exist $q\in\mathbb{Z}$ and $r\in\left\{  0,1,\ldots,N-1\right\}  $
such that $m-1=qN+r$ (namely, $q=q_{0}-1$ and $r=N-1$). In other words, the
statement $\mathcal{A}\left(  m-1\right)  $ holds\footnote{since the statement
$\mathcal{A}\left(  m-1\right)  $ is defined as \newline$\left(  \text{there
exist }q\in\mathbb{Z}\text{ and }r\in\left\{  0,1,\ldots,N-1\right\}  \text{
such that }m-1=qN+r\right)  $}. Thus, we have proven that $\mathcal{A}\left(
m-1\right)  $ holds in Case 1.

Let us next consider Case 2. In this case, we have $r_{0}\neq0$. Combining
$r_{0}\in\left\{  0,1,\ldots,N-1\right\}  $ with $r_{0}\neq0$, we obtain%
\[
r_{0}\in\left\{  0,1,\ldots,N-1\right\}  \setminus\left\{  0\right\}
=\left\{  1,2,\ldots,N-1\right\}  ,
\]
so that $r_{0}-1\in\left\{  0,1,\ldots,N-2\right\}  \subseteq\left\{
0,1,\ldots,N-1\right\}  $. Also,%
\[
\underbrace{m}_{=q_{0}N+r_{0}}-1=q_{0}N+r_{0}-1=q_{0}N+\left(  r_{0}-1\right)
.
\]
Since $r_{0}-1\in\left\{  0,1,\ldots,N-1\right\}  $, we can therefore conclude
that there exist $q\in\mathbb{Z}$ and $r\in\left\{  0,1,\ldots,N-1\right\}  $
such that $m-1=qN+r$ (namely, $q=q_{0}$ and $r=r_{0}-1$). In other words, the
statement $\mathcal{A}\left(  m-1\right)  $ holds\footnote{since the statement
$\mathcal{A}\left(  m-1\right)  $ is defined as \newline$\left(  \text{there
exist }q\in\mathbb{Z}\text{ and }r\in\left\{  0,1,\ldots,N-1\right\}  \text{
such that }m-1=qN+r\right)  $}. Thus, we have proven that $\mathcal{A}\left(
m-1\right)  $ holds in Case 2.

We have now proven that $\mathcal{A}\left(  m-1\right)  $ holds in each of the
two Cases 1 and 2. Since these two Cases cover all possibilities, we thus
conclude that $\mathcal{A}\left(  m-1\right)  $ always holds.

Now, forget that we fixed $m$. We thus have shown that if $m\in\mathbb{Z}%
_{\leq g}$ is such that $\mathcal{A}\left(  m\right)  $ holds, then
$\mathcal{A}\left(  m-1\right)  $ also holds. In other words, Assumption 3 is satisfied.]

We have now shown that all three Assumptions 1, 2 and 3 of Theorem
\ref{thm.ind.IPg+-} are satisfied. Thus, Theorem \ref{thm.ind.IPg+-} yields
that%
\begin{equation}
\mathcal{A}\left(  n\right)  \text{ holds for each }n\in\mathbb{Z}.
\label{pf.prop.ind.quo-rem-ex.at}%
\end{equation}

Now, let $n\in\mathbb{Z}$. Then, $\mathcal{A}\left(  n\right)  $ holds (by
(\ref{pf.prop.ind.quo-rem-ex.at})). In other words, there exist $q\in
\mathbb{Z}$ and $r\in\left\{  0,1,\ldots,N-1\right\}  $ such that $n=qN+r$
(because the statement $\mathcal{A}\left(  n\right)  $ is defined as $\left(
\text{there exist }q\in\mathbb{Z}\text{ and }r\in\left\{  0,1,\ldots
,N-1\right\}  \text{ such that }n=qN+r\right)  $). This proves Proposition
\ref{prop.ind.quo-rem-ex}.
\end{proof}

\subsubsection{Division with remainder}

We need one more lemma:

\begin{lemma}
\label{lem.ind.quo-rem-uni}Let $N$ be a positive integer. Let $\left(
q_{1},r_{1}\right)  \in\mathbb{Z}\times\left\{  0,1,\ldots,N-1\right\}  $ and
$\left(  q_{2},r_{2}\right)  \in\mathbb{Z}\times\left\{  0,1,\ldots
,N-1\right\}  $. Assume that $q_{1}N+r_{1}=q_{2}N+r_{2}$. Then, $\left(
q_{1},r_{1}\right)  =\left(  q_{2},r_{2}\right)  $.
\end{lemma}

\begin{proof}
[Proof of Lemma \ref{lem.ind.quo-rem-uni}.]We have $\left(  q_{1}%
,r_{1}\right)  \in\mathbb{Z}\times\left\{  0,1,\ldots,N-1\right\}  $. Thus,
$q_{1}\in\mathbb{Z}$ and $r_{1}\in\left\{  0,1,\ldots,N-1\right\}  $.
Similarly, $q_{2}\in\mathbb{Z}$ and $r_{2}\in\left\{  0,1,\ldots,N-1\right\}
$.

From $r_{1}\in\left\{  0,1,\ldots,N-1\right\}  $, we obtain $r_{1}\geq0$. From
$r_{2}\in\left\{  0,1,\ldots,N-1\right\}  $, we obtain $r_{2}\leq N-1$. Hence,
$\underbrace{r_{2}}_{\leq N-1}-\underbrace{r_{1}}_{\geq0}\leq\left(
N-1\right)  -0=N-1<N$.

Next, we shall prove that
\begin{equation}
q_{1}\leq q_{2}. \label{pf.lem.ind.quo-rem-uni.leq}%
\end{equation}

[\textit{Proof of (\ref{pf.lem.ind.quo-rem-uni.leq}):} Assume the contrary.
Thus, $q_{1}>q_{2}$, so that $q_{1}-q_{2}>0$. Hence, $q_{1}-q_{2}\geq1$ (since
$q_{1}-q_{2}$ is an integer). Therefore, $q_{1}-q_{2}-1\geq0$, so that
$N\left(  q_{1}-q_{2}-1\right)  \geq0$ (because $N>0$ and $q_{1}-q_{2}-1\geq0$).

But $q_{1}N+r_{1}=q_{2}N+r_{2}$, so that $q_{2}N+r_{2}=q_{1}N+r_{1}$. Hence,%
\[
r_{2}-r_{1}=q_{1}N-q_{2}N=N\left(  q_{1}-q_{2}\right)  =\underbrace{N\left(
q_{1}-q_{2}-1\right)  }_{\geq0}+N\geq N.
\]
This contradicts $r_{2}-r_{1}<N$. This contradiction shows that our assumption
was wrong. Hence, (\ref{pf.lem.ind.quo-rem-uni.leq}) is proven.]

Thus, we have proven that $q_{1}\leq q_{2}$. The same argument (with the roles
of $\left(  q_{1},r_{1}\right)  $ and $\left(  q_{2},r_{2}\right)  $
interchanged) shows that $q_{2}\leq q_{1}$. Combining the inequalities
$q_{1}\leq q_{2}$ and $q_{2}\leq q_{1}$, we obtain $q_{1}=q_{2}$.

Also, $q_{2}N+r_{2}=\underbrace{q_{1}}_{=q_{2}}N+r_{1}=q_{2}N+r_{1}$.
Subtracting $q_{2}N$ from both sides of this equality, we find $r_{2}=r_{1}$.
Hence, $r_{1}=r_{2}$.

Thus, $\left(  \underbrace{q_{1}}_{=q_{2}},\underbrace{r_{1}}_{=r_{2}}\right)
=\left(  q_{2},r_{2}\right)  $. This proves Lemma \ref{lem.ind.quo-rem-uni}.
\end{proof}

As we have already mentioned, Proposition \ref{prop.ind.quo-rem-ex} is just a
part of a crucial result from number theory:

\begin{theorem}
\label{thm.ind.quo-rem}Let $N$ be a positive integer. Let $n\in\mathbb{Z}$.
Then, there is a unique pair $\left(  q,r\right)  \in\mathbb{Z}\times\left\{
0,1,\ldots,N-1\right\}  $ such that $n=qN+r$.
\end{theorem}

Proving this theorem will turn out rather easy, since we have already done the
hard work with our proofs of Proposition \ref{prop.ind.quo-rem-ex} and Lemma
\ref{lem.ind.quo-rem-uni}:

\begin{noncompile}
(The following has been removed since I didn't end up using it.) We need one
more lemma:

\begin{lemma}
\label{lem.ind.divi-geq}Let $a$ and $d$ be two integers such that $a>0$ and
$d\mid a$. Then, $a\geq d$.
\end{lemma}

\begin{proof}
[Proof of Lemma \ref{lem.ind.divi-geq}.]We must prove that $a\geq d$. If
$d\leq0$, then this is obvious (because if $d\leq0$, then $a>0\geq d$). Hence,
for the rest of this proof, we WLOG assume that we don't have $d\leq0$. Thus,
we have $d>0$.

We have $d\mid a$. In other words, there exists an integer $w$ such that
$a=dw$ (by the definition of \textquotedblleft divides\textquotedblright).
Consider this $w$. If we had $w\leq0$, then we would have $dw\leq0$ (because
$d>0$ and $w\leq0$), which would contradict $dw=a>0$. Thus, we cannot have
$w\leq0$. Hence, we have $w>0$. Thus, $w\geq1$ (since $w$ is an integer). In
other words, $w-1\geq0$. From $d>0$ and $w-1\geq0$, we obtain $d\left(
w-1\right)  \geq0$. Now, $a=dw=\underbrace{d\left(  w-1\right)  }_{\geq
0}+d\geq d$. This proves Lemma \ref{lem.ind.divi-geq}.
\end{proof}
\end{noncompile}

\begin{proof}
[Proof of Theorem \ref{thm.ind.quo-rem}.]Proposition \ref{prop.ind.quo-rem-ex}
shows that there exist $q\in\mathbb{Z}$ and $r\in\left\{  0,1,\ldots
,N-1\right\}  $ such that $n=qN+r$. Consider these $q$ and $r$, and denote
them by $q_{0}$ and $r_{0}$. Thus, $q_{0}\in\mathbb{Z}$ and $r_{0}\in\left\{
0,1,\ldots,N-1\right\}  $ and $n=q_{0}N+r_{0}$. From $q_{0}\in\mathbb{Z}$ and
$r_{0}\in\left\{  0,1,\ldots,N-1\right\}  $, we obtain $\left(  q_{0}%
,r_{0}\right)  \in\mathbb{Z}\times\left\{  0,1,\ldots,N-1\right\}  $. Hence,
there exists \textbf{at least one} pair $\left(  q,r\right)  \in
\mathbb{Z}\times\left\{  0,1,\ldots,N-1\right\}  $ such that $n=qN+r$ (namely,
$\left(  q,r\right)  =\left(  q_{0},r_{0}\right)  $).

Now, let $\left(  q_{1},r_{1}\right)  $ and $\left(  q_{2},r_{2}\right)  $ be
two pairs $\left(  q,r\right)  \in\mathbb{Z}\times\left\{  0,1,\ldots
,N-1\right\}  $ such that $n=qN+r$. We shall prove that $\left(  q_{1}%
,r_{1}\right)  =\left(  q_{2},r_{2}\right)  $.

We have assumed that $\left(  q_{1},r_{1}\right)  $ is a pair $\left(
q,r\right)  \in\mathbb{Z}\times\left\{  0,1,\ldots,N-1\right\}  $ such that
$n=qN+r$. In other words, $\left(  q_{1},r_{1}\right)  $ is a pair in
$\mathbb{Z}\times\left\{  0,1,\ldots,N-1\right\}  $ and satisfies
$n=q_{1}N+r_{1}$. Similarly, $\left(  q_{2},r_{2}\right)  $ is a pair in
$\mathbb{Z}\times\left\{  0,1,\ldots,N-1\right\}  $ and satisfies
$n=q_{2}N+r_{2}$.

Hence, $q_{1}N+r_{1}=n=q_{2}N+r_{2}$. Thus, Lemma \ref{lem.ind.quo-rem-uni}
yields $\left(  q_{1},r_{1}\right)  =\left(  q_{2},r_{2}\right)  $.

Let us now forget that we fixed $\left(  q_{1},r_{1}\right)  $ and $\left(
q_{2},r_{2}\right)  $. We thus have shown that if $\left(  q_{1},r_{1}\right)
$ and $\left(  q_{2},r_{2}\right)  $ are two pairs $\left(  q,r\right)
\in\mathbb{Z}\times\left\{  0,1,\ldots,N-1\right\}  $ such that $n=qN+r$, then
$\left(  q_{1},r_{1}\right)  =\left(  q_{2},r_{2}\right)  $. In other words,
any two pairs $\left(  q,r\right)  \in\mathbb{Z}\times\left\{  0,1,\ldots
,N-1\right\}  $ such that $n=qN+r$ must be equal. In other words, there exists
\textbf{at most one} pair $\left(  q,r\right)  \in\mathbb{Z}\times\left\{
0,1,\ldots,N-1\right\}  $ such that $n=qN+r$. Since we also know that there
exists \textbf{at least one} such pair, we can therefore conclude that there
exists \textbf{exactly one} such pair. In other words, there is a unique pair
$\left(  q,r\right)  \in\mathbb{Z}\times\left\{  0,1,\ldots,N-1\right\}  $
such that $n=qN+r$. This proves Theorem \ref{thm.ind.quo-rem}.
\end{proof}

\begin{definition}
Let $N$ be a positive integer. Let $n\in\mathbb{Z}$. Theorem
\ref{thm.ind.quo-rem} says that there is a unique pair $\left(  q,r\right)
\in\mathbb{Z}\times\left\{  0,1,\ldots,N-1\right\}  $ such that $n=qN+r$.
Consider this pair $\left(  q,r\right)  $. Then, $q$ is called the
\textit{quotient of the division of }$n$ \textit{by }$N$ (or the
\textit{quotient obtained when }$n$ \textit{is divided by }$N$), whereas $r$
is called the \textit{remainder of the division of }$n$ \textit{by }$N$ (or
the \textit{remainder obtained when }$n$ \textit{is divided by }$N$).
\end{definition}

For example, the quotient of the division of $7$ by $3$ is $2$, whereas the
remainder of the division of $7$ by $3$ is $1$ (because $\left(  2,1\right)  $
is a pair in $\mathbb{Z}\times\left\{  0,1,2\right\}  $ such that
$7=2\cdot3+1$).

We collect some basic properties of remainders:

\begin{corollary}
\label{cor.ind.quo-rem.remmod}Let $N$ be a positive integer. Let
$n\in\mathbb{Z}$. Let $n\%N$ denote the remainder of the division of $n$ by
$N$.

\textbf{(a)} Then, $n\%N\in\left\{  0,1,\ldots,N-1\right\}  $ and $n\%N\equiv
n\operatorname{mod}N$.

\textbf{(b)} We have $N\mid n$ if and only if $n\%N=0$.

\textbf{(c)} Let $c\in\left\{  0,1,\ldots,N-1\right\}  $ be such that $c\equiv
n\operatorname{mod}N$. Then, $c=n\%N$.
\end{corollary}

\begin{proof}
[Proof of Corollary \ref{cor.ind.quo-rem.remmod}.]Theorem
\ref{thm.ind.quo-rem} says that there is a unique pair $\left(  q,r\right)
\in\mathbb{Z}\times\left\{  0,1,\ldots,N-1\right\}  $ such that $n=qN+r$.
Consider this pair $\left(  q,r\right)  $. Then, the remainder of the division
of $n$ by $N$ is $r$ (because this is how this remainder was defined). In
other words, $n\%N$ is $r$ (since $n\%N$ is the remainder of the division of
$n$ by $N$). Thus, $n\%N=r$. But $N\mid qN$ (since $q$ is an integer), so that
$qN\equiv0\operatorname{mod}N$. Hence, $\underbrace{qN}_{\equiv
0\operatorname{mod}N}+r\equiv0+r=r\operatorname{mod}N$. Hence, $r\equiv
qN+r=n\operatorname{mod}N$, so that $n\%N=r\equiv n\operatorname{mod}N$.
Furthermore, $n\%N=r\in\left\{  0,1,\ldots,N-1\right\}  $ (since $\left(
q,r\right)  \in\mathbb{Z}\times\left\{  0,1,\ldots,N-1\right\}  $). This
completes the proof of Corollary \ref{cor.ind.quo-rem.remmod} \textbf{(a)}.

\textbf{(b)} We have the following implication:%
\begin{equation}
\left(  N\mid n\right)  \Longrightarrow\left(  n\%N=0\right)  .
\label{pf.cor.ind.quo-rem.remmod.1}%
\end{equation}

[\textit{Proof of (\ref{pf.cor.ind.quo-rem.remmod.1}):} Assume that $N\mid n$.
We must prove that $n\%N=0$.

We have $N\mid n$. In other words, there exists some integer $w$ such that
$n=Nw$. Consider this $w$.

We have $N-1\in\mathbb{N}$ (since $N$ is a positive integer), thus
$0\in\left\{  0,1,\ldots,N-1\right\}  $. From $w\in\mathbb{Z}$ and
$0\in\left\{  0,1,\ldots,N-1\right\}  $, we obtain $\left(  w,0\right)
\in\mathbb{Z}\times\left\{  0,1,\ldots,N-1\right\}  $. Also,
$wN+0=wN=Nw=n=qN+r$. Hence, Lemma \ref{lem.ind.quo-rem-uni} (applied to
$\left(  q_{1},r_{1}\right)  =\left(  w,0\right)  $ and $\left(  q_{2}%
,r_{2}\right)  =\left(  q,r\right)  $) yields $\left(  w,0\right)  =\left(
q,r\right)  $. In other words, $w=q$ and $0=r$. Hence, $r=0$, so that
$n\%N=r=0$. This proves the implication (\ref{pf.cor.ind.quo-rem.remmod.1}).]

Next, we have the following implication:%
\begin{equation}
\left(  n\%N=0\right)  \Longrightarrow\left(  N\mid n\right)  .
\label{pf.cor.ind.quo-rem.remmod.2}%
\end{equation}

[\textit{Proof of (\ref{pf.cor.ind.quo-rem.remmod.2}):} Assume that $n\%N=0$.
We must prove that $N\mid n$.

We have $n=qN+\underbrace{r}_{=n\%N=0}=qN$. Thus, $N\mid n$. This proves the
implication (\ref{pf.cor.ind.quo-rem.remmod.2}).]

Combining the two implications (\ref{pf.cor.ind.quo-rem.remmod.1}) and
(\ref{pf.cor.ind.quo-rem.remmod.2}), we obtain the logical equivalence
$\left(  N\mid n\right)  \Longleftrightarrow\left(  n\%N=0\right)  $. In other
words, we have $N\mid n$ if and only if $n\%N=0$. This proves Corollary
\ref{cor.ind.quo-rem.remmod} \textbf{(b)}.

\textbf{(c)} We have $c\equiv n\operatorname{mod}N$. In other words, $N\mid
c-n$. In other words, there exists some integer $w$ such that $c-n=Nw$.
Consider this $w$.

From $-w\in\mathbb{Z}$ and $c\in\left\{  0,1,\ldots,N-1\right\}  $, we obtain
$\left(  -w,c\right)  \in\mathbb{Z}\times\left\{  0,1,\ldots,N-1\right\}  $.
Also, from $c-n=Nw$, we obtain $n=c-Nw=\left(  -w\right)  N+c$, so that
$\left(  -w\right)  N+c=n=qN+r$. Hence, Lemma \ref{lem.ind.quo-rem-uni}
(applied to $\left(  q_{1},r_{1}\right)  =\left(  -w,c\right)  $ and $\left(
q_{2},r_{2}\right)  =\left(  q,r\right)  $) yields $\left(  -w,c\right)
=\left(  q,r\right)  $. In other words, $-w=q$ and $c=r$. Hence, $c=r=n\%N$.
This proves Corollary \ref{cor.ind.quo-rem.remmod} \textbf{(c)}.
\end{proof}

Note that parts \textbf{(a)} and \textbf{(c)} of Corollary
\ref{cor.ind.quo-rem.remmod} (taken together) characterize the remainder
$n\%N$ as the unique element of $\left\{  0,1,\ldots,N-1\right\}  $ that is
congruent to $n$ modulo $N$. Corollary \ref{cor.ind.quo-rem.remmod}
\textbf{(b)} provides a simple algorithm to check whether a given integer $n$
is divisible by a given positive integer $N$; namely, it suffices to compute
the remainder $n\%N$ and check whether $n\%N=0$.

Let us further illustrate the usefulness of Theorem \ref{thm.ind.quo-rem} by
proving a fundamental property of odd numbers. Recall the following standard definitions:

\begin{definition}
Let $n\in\mathbb{Z}$.

\textbf{(a)} We say that the integer $n$ is \textit{even} if and only if $n$
is divisible by $2$.

\textbf{(b)} We say that the integer $n$ is \textit{odd} if and only if $n$ is
not divisible by $2$.
\end{definition}

This definition shows that any integer $n$ is either even or odd (but not both
at the same time).

It is clear that an integer $n$ is even if and only if it can be written in
the form $n=2m$ for some $m\in\mathbb{Z}$. Moreover, this $m$ is unique
(because $n=2m$ implies $m=n/2$). Let us prove a similar property for odd numbers:

\begin{proposition}
\label{prop.ind.quo-rem.odd}Let $n\in\mathbb{Z}$.

\textbf{(a)} The integer $n$ is odd if and only if $n$ can be written in the
form $n=2m+1$ for some $m\in\mathbb{Z}$.

\textbf{(b)} This $m$ is unique if it exists. (That is, any two integers
$m\in\mathbb{Z}$ satisfying $n=2m+1$ must be equal.)
\end{proposition}

We shall use Theorem \ref{thm.ind.quo-rem} several times in the below proof
(far more than necessary), mostly to illustrate how it can be applied.

\begin{proof}
[Proof of Proposition \ref{prop.ind.quo-rem.odd}.]\textbf{(a)} Let us first
prove the logical implication%
\begin{equation}
\left(  n\text{ is odd}\right)  \ \Longrightarrow\ \left(  \text{there exists
an }m\in\mathbb{Z}\text{ such that }n=2m+1\right)  .
\label{pf.prop.ind.quo-rem.odd.a.1}%
\end{equation}

[\textit{Proof of (\ref{pf.prop.ind.quo-rem.odd.a.1}):} Assume that $n$ is
odd. We must prove that there exists an $m\in\mathbb{Z}$ such that $n=2m+1$.

Theorem \ref{thm.ind.quo-rem} (applied to $N=2$) yields that there is a unique
pair $\left(  q,r\right)  \in\mathbb{Z}\times\left\{  0,1,\ldots,2-1\right\}
$ such that $n=q\cdot2+r$. Consider this $\left(  q,r\right)  $. From $\left(
q,r\right)  \in\mathbb{Z}\times\left\{  0,1,\ldots,2-1\right\}  $, we obtain
$q\in\mathbb{Z}$ and $r\in\left\{  0,1,\ldots,2-1\right\}  =\left\{
0,1\right\}  $.

We know that $n$ is odd; in other words, $n$ is not divisible by $2$ (by the
definition of \textquotedblleft odd\textquotedblright). If we had $n=2q$, then
$n$ would be divisible by $2$, which would contradict the fact that $n$ is not
divisible by $2$. Hence, we cannot have $n=2q$. If we had $r=0$, then we would
have $n=\underbrace{q\cdot2}_{=2q}+\underbrace{r}_{=0}=2q$, which would
contradict the fact that we cannot have $n=2q$. Hence, we cannot have $r=0$.
Thus, $r\neq0$.

Combining $r\in\left\{  0,1\right\}  $ with $r\neq0$, we obtain $r\in\left\{
0,1\right\}  \setminus\left\{  0\right\}  =\left\{  1\right\}  $. Thus, $r=1$.
Hence, $n=\underbrace{q\cdot2}_{=2q}+\underbrace{r}_{=1}=2q+1$. Thus, there
exists an $m\in\mathbb{Z}$ such that $n=2m+1$ (namely, $m=q$). This proves the
implication (\ref{pf.prop.ind.quo-rem.odd.a.1}).]

Next, we shall prove the logical implication%
\begin{equation}
\left(  \text{there exists an }m\in\mathbb{Z}\text{ such that }n=2m+1\right)
\ \Longrightarrow\ \left(  n\text{ is odd}\right)  .
\label{pf.prop.ind.quo-rem.odd.a.2}%
\end{equation}

[\textit{Proof of (\ref{pf.prop.ind.quo-rem.odd.a.2}):} Assume that there
exists an $m\in\mathbb{Z}$ such that $n=2m+1$. We must prove that $n$ is odd.

We have assumed that there exists an $m\in\mathbb{Z}$ such that $n=2m+1$.
Consider this $m$. Thus, the pair $\left(  m,1\right)  $ belongs to
$\mathbb{Z}\times\left\{  0,1,\ldots,2-1\right\}  $ (since $m\in\mathbb{Z}$
and $1\in\left\{  0,1,\ldots,2-1\right\}  $) and satisfies $n=m\cdot2+1$
(since $n=\underbrace{2m}_{=m\cdot2}+1=m\cdot2+1$). In other words, the pair
$\left(  m,1\right)  $ is a pair $\left(  q,r\right)  \in\mathbb{Z}%
\times\left\{  0,1,\ldots,2-1\right\}  $ such that $n=q\cdot2+r$.

Now, assume (for the sake of contradiction) that $n$ is divisible by $2$.
Thus, there exists some integer $w$ such that $n=2w$. Consider this $w$. Thus,
the pair $\left(  w,0\right)  $ belongs to $\mathbb{Z}\times\left\{
0,1,\ldots,2-1\right\}  $ (since $w\in\mathbb{Z}$ and $0\in\left\{
0,1,\ldots,2-1\right\}  $) and satisfies $n=w\cdot2+0$ (since $n=2w=w\cdot
2=w\cdot2+0$). In other words, the pair $\left(  w,0\right)  $ is a pair
$\left(  q,r\right)  \in\mathbb{Z}\times\left\{  0,1,\ldots,2-1\right\}  $
such that $n=q\cdot2+r$.

Theorem \ref{thm.ind.quo-rem} (applied to $N=2$) yields that there is a unique
pair $\left(  q,r\right)  \in\mathbb{Z}\times\left\{  0,1,\ldots,2-1\right\}
$ such that $n=q\cdot2+r$. Thus, there exists \textbf{at most} one such pair.
In other words, any two such pairs must be equal. Hence, the two pairs
$\left(  m,1\right)  $ and $\left(  w,0\right)  $ must be equal (since
$\left(  m,1\right)  $ and $\left(  w,0\right)  $ are two pairs $\left(
q,r\right)  \in\mathbb{Z}\times\left\{  0,1,\ldots,2-1\right\}  $ such that
$n=q\cdot2+r$). In other words, $\left(  m,1\right)  =\left(  w,0\right)  $.
In other words, $m=w$ and $1=0$. But $1=0$ is clearly absurd. Thus, we have
obtained a contradiction. This shows that our assumption (that $n$ is
divisible by $2$) was wrong. Hence, $n$ is not divisible by $2$. In other
words, $n$ is odd (by the definition of \textquotedblleft
odd\textquotedblright). This proves the implication
(\ref{pf.prop.ind.quo-rem.odd.a.2}).]

Combining the two implications (\ref{pf.prop.ind.quo-rem.odd.a.1}) and
(\ref{pf.prop.ind.quo-rem.odd.a.2}), we obtain the logical equivalence%
\begin{align*}
\left(  n\text{ is odd}\right)  \  &  \Longleftrightarrow\ \left(  \text{there
exists an }m\in\mathbb{Z}\text{ such that }n=2m+1\right) \\
&  \Longleftrightarrow\ \left(  n\text{ can be written in the form
}n=2m+1\text{ for some }m\in\mathbb{Z}\right)  .
\end{align*}
In other words, the integer $n$ is odd if and only if $n$ can be written in
the form $n=2m+1$ for some $m\in\mathbb{Z}$. This proves Proposition
\ref{prop.ind.quo-rem.odd} \textbf{(a)}.

\textbf{(b)} This is easy to prove in any way, but let us prove this using
Theorem \ref{thm.ind.quo-rem} just in order to illustrate the use of the
latter theorem.

We must prove that any two integers $m\in\mathbb{Z}$ satisfying $n=2m+1$ must
be equal.

Let $m_{1}$ and $m_{2}$ be two integers $m\in\mathbb{Z}$ satisfying $n=2m+1$.
We shall show that $m_{1}=m_{2}$.

We know that $m_{1}$ is an integer $m\in\mathbb{Z}$ satisfying $n=2m+1$. In
other words, $m_{1}$ is an integer in $\mathbb{Z}$ and satisfies $n=2m_{1}+1$.
Thus, the pair $\left(  m_{1},1\right)  $ belongs to $\mathbb{Z}\times\left\{
0,1,\ldots,2-1\right\}  $ (since $m_{1}\in\mathbb{Z}$ and $1\in\left\{
0,1,\ldots,2-1\right\}  $) and satisfies $n=m_{1}\cdot2+1$ (since
$n=\underbrace{2m_{1}}_{=m_{1}\cdot2}+1=m_{1}\cdot2+1$). In other words, the
pair $\left(  m_{1},1\right)  $ is a pair $\left(  q,r\right)  \in
\mathbb{Z}\times\left\{  0,1,\ldots,2-1\right\}  $ such that $n=q\cdot2+r$.
The same argument (applied to $m_{2}$ instead of $m_{1}$) shows that $\left(
m_{2},1\right)  $ is a pair $\left(  q,r\right)  \in\mathbb{Z}\times\left\{
0,1,\ldots,2-1\right\}  $ such that $n=q\cdot2+r$.

Theorem \ref{thm.ind.quo-rem} (applied to $N=2$) yields that there is a unique
pair $\left(  q,r\right)  \in\mathbb{Z}\times\left\{  0,1,\ldots,2-1\right\}
$ such that $n=q\cdot2+r$. Thus, there exists \textbf{at most} one such pair.
In other words, any two such pairs must be equal. Hence, the two pairs
$\left(  m_{1},1\right)  $ and $\left(  m_{2},1\right)  $ must be equal (since
$\left(  m_{1},1\right)  $ and $\left(  m_{2},1\right)  $ are two pairs
$\left(  q,r\right)  \in\mathbb{Z}\times\left\{  0,1,\ldots,2-1\right\}  $
such that $n=q\cdot2+r$). In other words, $\left(  m_{1},1\right)  =\left(
m_{2},1\right)  $. In other words, $m_{1}=m_{2}$ and $1=1$. Hence, we have
shown that $m_{1}=m_{2}$.

Now, forget that we fixed $m_{1}$ and $m_{2}$. We thus have proven that if
$m_{1}$ and $m_{2}$ are two integers $m\in\mathbb{Z}$ satisfying $n=2m+1$,
then $m_{1}=m_{2}$. In other words, any two integers $m\in\mathbb{Z}$
satisfying $n=2m+1$ must be equal. In other words, the $m$ in Proposition
\ref{prop.ind.quo-rem.odd} \textbf{(a)} is unique. This proves Proposition
\ref{prop.ind.quo-rem.odd} \textbf{(b)}.
\end{proof}

We can use this to obtain the following fundamental fact:

\begin{corollary}
\label{cor.mod.-1powers}Let $n\in\mathbb{Z}$.

\textbf{(a)} If $n$ is even, then $\left(  -1\right)  ^{n}=1$.

\textbf{(b)} If $n$ is odd, then $\left(  -1\right)  ^{n}=-1$.
\end{corollary}

\begin{proof}
[Proof of Corollary \ref{cor.mod.-1powers}.]\textbf{(a)} Assume that $n$ is
even. In other words, $n$ is divisible by $2$ (by the definition of
\textquotedblleft even\textquotedblright). In other words, $2\mid n$. In other
words, there exists an integer $w$ such that $n=2w$. Consider this $w$. From
$n=2w$, we obtain $\left(  -1\right)  ^{n}=\left(  -1\right)  ^{2w}=\left(
\underbrace{\left(  -1\right)  ^{2}}_{=1}\right)  ^{w}=1^{w}=1$. This proves
Corollary \ref{cor.mod.-1powers} \textbf{(a)}.

\textbf{(b)} Assume that $n$ is odd. Proposition \ref{prop.ind.quo-rem.odd}
\textbf{(a)} shows that the integer $n$ is odd if and only if $n$ can be
written in the form $n=2m+1$ for some $m\in\mathbb{Z}$. Hence, $n$ can be
written in the form $n=2m+1$ for some $m\in\mathbb{Z}$ (since the integer $n$
is odd). Consider this $m$. From $n=2m+1$, we obtain%
\[
\left(  -1\right)  ^{n}=\left(  -1\right)  ^{2m+1}=\left(  -1\right)
^{2m}\left(  -1\right)  =-\underbrace{\left(  -1\right)  ^{2m}}_{=\left(
\left(  -1\right)  ^{2}\right)  ^{m}}=-\left(  \underbrace{\left(  -1\right)
^{2}}_{=1}\right)  ^{m}=-\underbrace{1^{m}}_{=1}=-1.
\]
This proves Corollary \ref{cor.mod.-1powers} \textbf{(b)}.
\end{proof}

Let us state two more fundamental facts, which are proven in Exercise
\ref{exe.mod.parity}:

\begin{proposition}
\label{prop.mod.parity}Let $u$ and $v$ be two integers. Then, we have the
following chain of logical equivalences:%
\begin{align*}
\left(  u\equiv v\operatorname{mod}2\right)  \  &  \Longleftrightarrow
\ \left(  u\text{ and }v\text{ are either both even or both odd}\right) \\
&  \Longleftrightarrow\ \left(  \left(  -1\right)  ^{u}=\left(  -1\right)
^{v}\right)  .
\end{align*}

\end{proposition}

\begin{proposition}
\label{prop.mod.even-odd-mod2}Let $n\in\mathbb{Z}$.

\textbf{(a)} The integer $n$ is even if and only if $n\equiv
0\operatorname{mod}2$.

\textbf{(b)} The integer $n$ is odd if and only if $n\equiv1\operatorname{mod}%
2$.
\end{proposition}

\begin{exercise}
\label{exe.mod.parity}Prove Proposition \ref{prop.mod.parity} and Proposition
\ref{prop.mod.even-odd-mod2}.
\end{exercise}

\begin{exercise}
\label{exe.mod.unique-cong}Let $N$ be a positive integer. Let $p\in\mathbb{Z}$
and $h\in\mathbb{Z}$. Prove that there exists a \textbf{unique} element
$g\in\left\{  p+1,p+2,\ldots,p+N\right\}  $ satisfying $g\equiv
h\operatorname{mod}N$.
\end{exercise}

\begin{exercise}
\label{exe.mod.ak-bk}Let $k\in\mathbb{N}$, $a\in\mathbb{Z}$ and $b\in
\mathbb{Z}$.

\textbf{(a)} Prove that $a-b\mid a^{k}-b^{k}$.

\textbf{(b)} Assume that $k$ is odd. Prove that $a+b\mid a^{k}+b^{k}$.
\end{exercise}

\begin{exercise}
\label{exe.ind.LP2.another-div}Fix an \textbf{odd} positive integer $r$.
Consider the sequence $\left(  b_{0},b_{1},b_{2},\ldots\right)  $ defined in
Proposition \ref{prop.ind.LP2}.

Prove that $b_{n}+1\mid b_{n-1}\left(  b_{n+2}+1\right)  $ for each positive
integer $n$. (Note that the statement \textquotedblleft$b_{n}+1\mid
b_{n-1}\left(  b_{n+2}+1\right)  $\textquotedblright\ makes sense, since
Proposition \ref{prop.ind.LP2} \textbf{(a)} yields that all three numbers
$b_{n},b_{n-1},b_{n+2}$ belong to $\mathbb{N}$.)
\end{exercise}

\subsubsection{Backwards induction principles}

When we use Theorem \ref{thm.ind.IPg+-} to prove a statement, we can regard
the proof of Assumption 2 as a (regular) induction step (\textquotedblleft
forwards induction step\textquotedblright), and regard the proof of Assumption
3 as a sort of \textquotedblleft backwards induction step\textquotedblright.
There are also \textquotedblleft backwards induction
principles\textquotedblright\ which include a \textquotedblleft backwards
induction step\textquotedblright\ but no \textquotedblleft forwards induction
step\textquotedblright. Here are two such principles:

\begin{theorem}
\label{thm.ind.IPg-}Let $g\in\mathbb{Z}$. Let $\mathbb{Z}_{\leq g}$ be the set
$\left\{  g,g-1,g-2,\ldots\right\}  $ (that is, the set of all integers that
are $\leq g$). For each $n\in\mathbb{Z}_{\leq g}$, let $\mathcal{A}\left(
n\right)  $ be a logical statement.

Assume the following:

\begin{statement}
\textit{Assumption 1:} The statement $\mathcal{A}\left(  g\right)  $ holds.
\end{statement}

\begin{statement}
\textit{Assumption 2:} If $m\in\mathbb{Z}_{\leq g}$ is such that
$\mathcal{A}\left(  m\right)  $ holds, then $\mathcal{A}\left(  m-1\right)  $
also holds.
\end{statement}

Then, $\mathcal{A}\left(  n\right)  $ holds for each $n\in\mathbb{Z}_{\leq g}$.
\end{theorem}

\begin{theorem}
\label{thm.ind.IPgh-}Let $g\in\mathbb{Z}$ and $h\in\mathbb{Z}$. For each
$n\in\left\{  g,g+1,\ldots,h\right\}  $, let $\mathcal{A}\left(  n\right)  $
be a logical statement.

Assume the following:

\begin{statement}
\textit{Assumption 1:} If $g\leq h$, then the statement $\mathcal{A}\left(
h\right)  $ holds.
\end{statement}

\begin{statement}
\textit{Assumption 2:} If $m\in\left\{  g+1,g+2,\ldots,h\right\}  $ is such
that $\mathcal{A}\left(  m\right)  $ holds, then $\mathcal{A}\left(
m-1\right)  $ also holds.
\end{statement}

Then, $\mathcal{A}\left(  n\right)  $ holds for each $n\in\left\{
g,g+1,\ldots,h\right\}  $.
\end{theorem}

Theorem \ref{thm.ind.IPg-} is an analogue of Theorem \ref{thm.ind.IPg}, while
Theorem \ref{thm.ind.IPgh-} is an analogue of Theorem \ref{thm.ind.IPgh}.
However, it is not hard to derive these theorems from the induction principles
we already know:

\begin{exercise}
\label{exe.ind.backw}Prove Theorem \ref{thm.ind.IPg-} and Theorem
\ref{thm.ind.IPgh-}.
\end{exercise}

It is also easy to state and prove a \textquotedblleft
backwards\textquotedblright\ analogue of Theorem \ref{thm.ind.SIP}. (We leave
this to the reader.)

A proof using Theorem \ref{thm.ind.IPg-} or using Theorem \ref{thm.ind.IPgh-}
is usually called a \textit{proof by descending induction} or a \textit{proof
by backwards induction}.

\subsection{Induction from $k-1$ to $k$}

\subsubsection{The principle}

Let us next show yet another \textquotedblleft alternative induction
principle\textquotedblright, which differs from Theorem \ref{thm.ind.IPg} in a
mere notational detail:

\begin{theorem}
\label{thm.ind.IPg-1}Let $g\in\mathbb{Z}$. For each $n\in\mathbb{Z}_{\geq g}$,
let $\mathcal{A}\left(  n\right)  $ be a logical statement.

Assume the following:

\begin{statement}
\textit{Assumption 1:} The statement $\mathcal{A}\left(  g\right)  $ holds.
\end{statement}

\begin{statement}
\textit{Assumption 2:} If $k\in\mathbb{Z}_{\geq g+1}$ is such that
$\mathcal{A}\left(  k-1\right)  $ holds, then $\mathcal{A}\left(  k\right)  $
also holds.
\end{statement}

Then, $\mathcal{A}\left(  n\right)  $ holds for each $n\in\mathbb{Z}_{\geq g}$.
\end{theorem}

Roughly speaking, this Theorem \ref{thm.ind.IPg-1} is just Theorem
\ref{thm.ind.IPg}, except that the variable $m$ in Assumption 2 has been
renamed as $k-1$. Consequently, it stands to reason that Theorem
\ref{thm.ind.IPg-1} can easily be derived from Theorem \ref{thm.ind.IPg}. Here
is the derivation in full detail:

\begin{proof}
[Proof of Theorem \ref{thm.ind.IPg-1}.]For each $n\in\mathbb{Z}_{\geq g}$, we
define the logical statement $\mathcal{B}\left(  n\right)  $ to be the
statement $\mathcal{A}\left(  n\right)  $. Thus, $\mathcal{B}\left(  n\right)
=\mathcal{A}\left(  n\right)  $ for each $n\in\mathbb{Z}_{\geq g}$. Applying
this to $n=g$, we obtain $\mathcal{B}\left(  g\right)  =\mathcal{A}\left(
g\right)  $ (since $g\in\mathbb{Z}_{\geq g}$).

We shall now show that the two Assumptions A and B of Corollary
\ref{cor.ind.IPg.renamed} are satisfied.

Indeed, recall that Assumption 1 is satisfied. In other words, the statement
$\mathcal{A}\left(  g\right)  $ holds. In other words, the statement
$\mathcal{B}\left(  g\right)  $ holds (since $\mathcal{B}\left(  g\right)
=\mathcal{A}\left(  g\right)  $). In other words, Assumption A is satisfied.

We shall next show that Assumption B is satisfied. Indeed, let $p\in
\mathbb{Z}_{\geq g}$ be such that $\mathcal{B}\left(  p\right)  $ holds.
Recall that the statement $\mathcal{B}\left(  p\right)  $ was defined to be
the statement $\mathcal{A}\left(  p\right)  $. Thus, $\mathcal{B}\left(
p\right)  =\mathcal{A}\left(  p\right)  $. Hence, $\mathcal{A}\left(
p\right)  $ holds (since $\mathcal{B}\left(  p\right)  $ holds). Now, let
$k=p+1$. We know that $p\in\mathbb{Z}_{\geq g}$; in other words, $p$ is an
integer and satisfies $p\geq g$. Hence, $k=p+1$ is an integer as well and
satisfies $k=\underbrace{p}_{\geq g}+1\geq g+1$. In other words,
$k\in\mathbb{Z}_{\geq g+1}$. Moreover, from $k=p+1$, we obtain $k-1=p$. Hence,
$\mathcal{A}\left(  k-1\right)  =\mathcal{A}\left(  p\right)  $. Thus,
$\mathcal{A}\left(  k-1\right)  $ holds (since $\mathcal{A}\left(  p\right)  $
holds). Thus, Assumption 2 shows that $\mathcal{A}\left(  k\right)  $ also
holds. But the statement $\mathcal{B}\left(  k\right)  $ was defined to be the
statement $\mathcal{A}\left(  k\right)  $. Hence, $\mathcal{B}\left(
k\right)  =\mathcal{A}\left(  k\right)  $, so that $\mathcal{A}\left(
k\right)  =\mathcal{B}\left(  k\right)  =\mathcal{B}\left(  p+1\right)  $
(since $k=p+1$). Thus, the statement $\mathcal{B}\left(  p+1\right)  $ holds
(since $\mathcal{A}\left(  k\right)  $ holds). Now, forget that we fixed $p$.
We thus have shown that if $p\in\mathbb{Z}_{\geq g}$ is such that
$\mathcal{B}\left(  p\right)  $ holds, then $\mathcal{B}\left(  p+1\right)  $
also holds. In other words, Assumption B is satisfied.

We have now proven that both Assumptions A and B of Corollary
\ref{cor.ind.IPg.renamed} are satisfied. Hence, Corollary
\ref{cor.ind.IPg.renamed} shows that $\mathcal{B}\left(  n\right)  $ holds for
each $n\in\mathbb{Z}_{\geq g}$. In other words, $\mathcal{A}\left(  n\right)
$ holds for each $n\in\mathbb{Z}_{\geq g}$ (because each $n\in\mathbb{Z}_{\geq
g}$ satisfies $\mathcal{B}\left(  n\right)  =\mathcal{A}\left(  n\right)  $
(by the definition of $\mathcal{B}\left(  n\right)  $)). This proves Theorem
\ref{thm.ind.IPg-1}.
\end{proof}

Proofs that use Theorem \ref{thm.ind.IPg-1} are usually called \textit{proofs
by induction} or \textit{induction proofs}. As an example of such a proof, let
us show the following identity:

\begin{proposition}
\label{prop.ind.alt-harm}For every $n\in\mathbb{N}$, we have%
\begin{equation}
\sum_{i=1}^{2n}\dfrac{\left(  -1\right)  ^{i-1}}{i}=\sum_{i=n+1}^{2n}\dfrac
{1}{i}. \label{eq.prop.ind.alt-harm.claim}%
\end{equation}

\end{proposition}

The equality (\ref{eq.prop.ind.alt-harm.claim}) can be rewritten as%
\[
\dfrac{1}{1}-\dfrac{1}{2}+\dfrac{1}{3}-\dfrac{1}{4}\pm\cdots+\dfrac{1}%
{2n-1}-\dfrac{1}{2n}=\dfrac{1}{n+1}+\dfrac{1}{n+2}+\cdots+\dfrac{1}{2n}%
\]
(where all the signs on the right hand side are $+$ signs, whereas the signs
on the left hand side alternate between $+$ signs and $-$ signs).

\begin{proof}
[Proof of Proposition \ref{prop.ind.alt-harm}.]For each $n\in\mathbb{Z}%
_{\geq0}$, we let $\mathcal{A}\left(  n\right)  $ be the statement%
\[
\left(  \sum_{i=1}^{2n}\dfrac{\left(  -1\right)  ^{i-1}}{i}=\sum_{i=n+1}%
^{2n}\dfrac{1}{i}\right)  .
\]
Our next goal is to prove the statement $\mathcal{A}\left(  n\right)  $ for
each $n\in\mathbb{Z}_{\geq0}$.

We first notice that the statement $\mathcal{A}\left(  0\right)  $
holds\footnote{\textit{Proof.} We have $\sum_{i=1}^{2\cdot0}\dfrac{\left(
-1\right)  ^{i-1}}{i}=\left(  \text{empty sum}\right)  =0$. Comparing this
with $\sum_{i=0+1}^{2\cdot0}\dfrac{1}{i}=\left(  \text{empty sum}\right)  =0$,
we obtain $\sum_{i=1}^{2\cdot0}\dfrac{\left(  -1\right)  ^{i-1}}{i}%
=\sum_{i=0+1}^{2\cdot0}\dfrac{1}{i}$. But this is precisely the statement
$\mathcal{A}\left(  0\right)  $ (since $\mathcal{A}\left(  0\right)  $ is
defined to be the statement $\left(  \sum_{i=1}^{2\cdot0}\dfrac{\left(
-1\right)  ^{i-1}}{i}=\sum_{i=0+1}^{2\cdot0}\dfrac{1}{i}\right)  $). Hence,
the statement $\mathcal{A}\left(  0\right)  $ holds.}.

Now, we claim that
\begin{equation}
\text{if }k\in\mathbb{Z}_{\geq0+1}\text{ is such that }\mathcal{A}\left(
k-1\right)  \text{ holds, then }\mathcal{A}\left(  k\right)  \text{ also
holds.} \label{pf.prop.ind.alt-harm.step}%
\end{equation}

[\textit{Proof of (\ref{pf.prop.ind.alt-harm.step}):} Let $k\in\mathbb{Z}%
_{\geq0+1}$ be such that $\mathcal{A}\left(  k-1\right)  $ holds. We must show
that $\mathcal{A}\left(  k\right)  $ also holds.

We have $k\in\mathbb{Z}_{\geq0+1}$. Thus, $k$ is an integer and satisfies
$k\geq0+1=1$.

We have assumed that $\mathcal{A}\left(  k-1\right)  $ holds. In other words,
\begin{equation}
\sum_{i=1}^{2\left(  k-1\right)  }\dfrac{\left(  -1\right)  ^{i-1}}{i}%
=\sum_{i=\left(  k-1\right)  +1}^{2\left(  k-1\right)  }\dfrac{1}{i}
\label{pf.prop.ind.alt-harm.IH}%
\end{equation}
holds\footnote{because $\mathcal{A}\left(  k-1\right)  $ is defined to be the
statement $\left(  \sum_{i=1}^{2\left(  k-1\right)  }\dfrac{\left(  -1\right)
^{i-1}}{i}=\sum_{i=\left(  k-1\right)  +1}^{2\left(  k-1\right)  }\dfrac{1}%
{i}\right)  $}.

We have $\left(  -1\right)  ^{2\left(  k-1\right)  }=\left(
\underbrace{\left(  -1\right)  ^{2}}_{=1}\right)  ^{k-1}=1^{k-1}=1$. But
$2k-1=2\left(  k-1\right)  +1$. Thus, $\left(  -1\right)  ^{2k-1}=\left(
-1\right)  ^{2\left(  k-1\right)  +1}=\underbrace{\left(  -1\right)
^{2\left(  k-1\right)  }}_{=1}\underbrace{\left(  -1\right)  ^{1}}_{=-1}=-1$.

Now, $k\geq1$, so that $2k\geq2$ and therefore $2k-1\geq1$. Hence, we can
split off the addend for $i=2k-1$ from the sum $\sum_{i=1}^{2k-1}%
\dfrac{\left(  -1\right)  ^{i-1}}{i}$. We thus obtain%
\begin{align}
\sum_{i=1}^{2k-1}\dfrac{\left(  -1\right)  ^{i-1}}{i}  &  =\sum_{i=1}^{\left(
2k-1\right)  -1}\dfrac{\left(  -1\right)  ^{i-1}}{i}+\dfrac{\left(  -1\right)
^{\left(  2k-1\right)  -1}}{2k-1}\nonumber\\
&  =\underbrace{\sum_{i=1}^{2\left(  k-1\right)  }\dfrac{\left(  -1\right)
^{i-1}}{i}}_{\substack{=\sum_{i=\left(  k-1\right)  +1}^{2\left(  k-1\right)
}\dfrac{1}{i}\\\text{(by (\ref{pf.prop.ind.alt-harm.IH}))}}%
}+\underbrace{\dfrac{\left(  -1\right)  ^{2\left(  k-1\right)  }}{2k-1}%
}_{\substack{=\dfrac{1}{2k-1}\\\text{(since }\left(  -1\right)  ^{2\left(
k-1\right)  }=1\text{)}}}\nonumber\\
&  \ \ \ \ \ \ \ \ \ \ \left(  \text{since }\left(  2k-1\right)  -1=2\left(
k-1\right)  \right) \nonumber\\
&  =\sum_{i=\left(  k-1\right)  +1}^{2\left(  k-1\right)  }\dfrac{1}{i}%
+\dfrac{1}{2k-1}=\sum_{i=k}^{2k-2}\dfrac{1}{i}+\dfrac{1}{2k-1}
\label{pf.prop.ind.alt-harm.0}%
\end{align}
(since $\left(  k-1\right)  +1=k$ and $2\left(  k-1\right)  =2k-2$).

On the other hand, $2k\geq2\geq1$. Hence, we can split off the addend for
$i=2k$ from the sum $\sum_{i=1}^{2k}\dfrac{\left(  -1\right)  ^{i-1}}{i}$. We
thus obtain%
\begin{align}
\sum_{i=1}^{2k}\dfrac{\left(  -1\right)  ^{i-1}}{i}  &  =\underbrace{\sum
_{i=1}^{2k-1}\dfrac{\left(  -1\right)  ^{i-1}}{i}}_{\substack{=\sum
_{i=k}^{2k-2}\dfrac{1}{i}+\dfrac{1}{2k-1}\\\text{(by
(\ref{pf.prop.ind.alt-harm.0}))}}}+\underbrace{\dfrac{\left(  -1\right)
^{2k-1}}{2k}}_{\substack{=\dfrac{-1}{2k}\\\text{(since }\left(  -1\right)
^{2k-1}=-1\text{)}}}\nonumber\\
&  =\sum_{i=k}^{2k-2}\dfrac{1}{i}+\dfrac{1}{2k-1}+\dfrac{-1}{2k}.
\label{pf.prop.ind.alt-harm.1}%
\end{align}

But we have $\left(  2k-1\right)  -k=k-1\geq0$ (since $k\geq1$). Thus,
$2k-1\geq k$. Thus, we can split off the addend for $i=2k-1$ from the sum
$\sum_{i=k}^{2k-1}\dfrac{1}{i}$. We thus obtain%
\begin{equation}
\sum_{i=k}^{2k-1}\dfrac{1}{i}=\sum_{i=k}^{\left(  2k-1\right)  -1}\dfrac{1}%
{i}+\dfrac{1}{2k-1}=\sum_{i=k}^{2k-2}\dfrac{1}{i}+\dfrac{1}{2k-1}
\label{pf.prop.ind.alt-harm.2}%
\end{equation}
(since $\left(  2k-1\right)  -1=2k-2$). Hence, (\ref{pf.prop.ind.alt-harm.1})
becomes%
\begin{equation}
\sum_{i=1}^{2k}\dfrac{\left(  -1\right)  ^{i-1}}{i}=\underbrace{\sum
_{i=k}^{2k-2}\dfrac{1}{i}+\dfrac{1}{2k-1}}_{\substack{=\sum_{i=k}^{2k-1}%
\dfrac{1}{i}\\\text{(by (\ref{pf.prop.ind.alt-harm.2}))}}}+\dfrac{-1}{2k}%
=\sum_{i=k}^{2k-1}\dfrac{1}{i}+\dfrac{-1}{2k}. \label{pf.prop.ind.alt-harm.3}%
\end{equation}

But we have $k+1\leq2k$ (since $2k-\left(  k+1\right)  =k-1\geq0$). Thus, we
can split off the addend for $i=2k$ from the sum $\sum_{i=k+1}^{2k}\dfrac
{1}{i}$. We thus obtain%
\[
\sum_{i=k+1}^{2k}\dfrac{1}{i}=\sum_{i=k+1}^{2k-1}\dfrac{1}{i}+\dfrac{1}{2k}.
\]
Hence,%
\begin{equation}
\sum_{i=k+1}^{2k-1}\dfrac{1}{i}=\sum_{i=k+1}^{2k}\dfrac{1}{i}-\dfrac{1}{2k}.
\label{pf.prop.ind.alt-harm.4}%
\end{equation}

Also, $k\leq2k-1$ (since $\left(  2k-1\right)  -k=k-1\geq0$). Thus, we can
split off the addend for $i=k$ from the sum $\sum_{i=k}^{2k-1}\dfrac{1}{i}$.
We thus obtain%
\[
\sum_{i=k}^{2k-1}\dfrac{1}{i}=\dfrac{1}{k}+\underbrace{\sum_{i=k+1}%
^{2k-1}\dfrac{1}{i}}_{\substack{=\sum_{i=k+1}^{2k}\dfrac{1}{i}-\dfrac{1}%
{2k}\\\text{(by (\ref{pf.prop.ind.alt-harm.4}))}}}=\dfrac{1}{k}+\sum
_{i=k+1}^{2k}\dfrac{1}{i}-\dfrac{1}{2k}=\sum_{i=k+1}^{2k}\dfrac{1}%
{i}+\underbrace{\dfrac{1}{k}-\dfrac{1}{2k}}_{=\dfrac{1}{2k}}=\sum_{i=k+1}%
^{2k}\dfrac{1}{i}+\dfrac{1}{2k}.
\]
Subtracting $\dfrac{1}{2k}$ from this equality, we obtain%
\[
\sum_{i=k}^{2k-1}\dfrac{1}{i}-\dfrac{1}{2k}=\sum_{i=k+1}^{2k}\dfrac{1}{i}.
\]
Hence,%
\[
\sum_{i=k+1}^{2k}\dfrac{1}{i}=\sum_{i=k}^{2k-1}\dfrac{1}{i}-\dfrac{1}{2k}%
=\sum_{i=k}^{2k-1}\dfrac{1}{i}+\dfrac{-1}{2k}.
\]
Comparing this with (\ref{pf.prop.ind.alt-harm.3}), we obtain%
\begin{equation}
\sum_{i=1}^{2k}\dfrac{\left(  -1\right)  ^{i-1}}{i}=\sum_{i=k+1}^{2k}\dfrac
{1}{i}. \label{pf.prop.ind.alt-harm.8}%
\end{equation}
But this is precisely the statement $\mathcal{A}\left(  k\right)
$\ \ \ \ \footnote{because $\mathcal{A}\left(  k\right)  $ is defined to be
the statement $\left(  \sum_{i=1}^{2k}\dfrac{\left(  -1\right)  ^{i-1}}%
{i}=\sum_{i=k+1}^{2k}\dfrac{1}{i}\right)  $}. Thus, the statement
$\mathcal{A}\left(  k\right)  $ holds.

Now, forget that we fixed $k$. We thus have shown that if $k\in\mathbb{Z}%
_{\geq0+1}$ is such that $\mathcal{A}\left(  k-1\right)  $ holds, then
$\mathcal{A}\left(  k\right)  $ also holds. This proves
(\ref{pf.prop.ind.alt-harm.step}).]

Now, both assumptions of Theorem \ref{thm.ind.IPg-1} (applied to $g=0$) are
satisfied (indeed, Assumption 1 holds because the statement $\mathcal{A}%
\left(  0\right)  $ holds, whereas Assumption 2 holds because of
(\ref{pf.prop.ind.alt-harm.step})). Thus, Theorem \ref{thm.ind.IPg-1} (applied
to $g=0$) shows that $\mathcal{A}\left(  n\right)  $ holds for each
$n\in\mathbb{Z}_{\geq0}$. In other words, $\sum_{i=1}^{2n}\dfrac{\left(
-1\right)  ^{i-1}}{i}=\sum_{i=n+1}^{2n}\dfrac{1}{i}$ holds for each
$n\in\mathbb{Z}_{\geq0}$ (since $\mathcal{A}\left(  n\right)  $ is the
statement $\left(  \sum_{i=1}^{2n}\dfrac{\left(  -1\right)  ^{i-1}}{i}%
=\sum_{i=n+1}^{2n}\dfrac{1}{i}\right)  $). In other words, $\sum_{i=1}%
^{2n}\dfrac{\left(  -1\right)  ^{i-1}}{i}=\sum_{i=n+1}^{2n}\dfrac{1}{i}$ holds
for each $n\in\mathbb{N}$ (because $\mathbb{Z}_{\geq0}=\mathbb{N}$). This
proves Proposition \ref{prop.ind.alt-harm}.
\end{proof}

\subsubsection{Conventions for writing proofs using \textquotedblleft$k-1$ to
$k$\textquotedblright\ induction}

Just like most of the other induction principles that we have so far
introduced, Theorem \ref{thm.ind.IPg-1} is not usually invoked explicitly when
it is used; instead, its use is signalled by certain words:

\begin{convention}
\label{conv.ind.IPg-1lang}Let $g\in\mathbb{Z}$. For each $n\in\mathbb{Z}_{\geq
g}$, let $\mathcal{A}\left(  n\right)  $ be a logical statement. Assume that
you want to prove that $\mathcal{A}\left(  n\right)  $ holds for each
$n\in\mathbb{Z}_{\geq g}$.

Theorem \ref{thm.ind.IPg-1} offers the following strategy for proving this:
First show that Assumption 1 of Theorem \ref{thm.ind.IPg-1} is satisfied;
then, show that Assumption 2 of Theorem \ref{thm.ind.IPg-1} is satisfied;
then, Theorem \ref{thm.ind.IPg-1} automatically completes your proof.

A proof that follows this strategy is called a \textit{proof by induction on
}$n$ (or \textit{proof by induction over }$n$) \textit{starting at }$g$ or
(less precisely) an \textit{inductive proof}. Most of the time, the words
\textquotedblleft starting at $g$\textquotedblright\ are omitted, since the
value of $g$ is usually clear from the statement that is being proven.
Usually, the statements $\mathcal{A}\left(  n\right)  $ are not explicitly
stated in the proof either, since they can also be inferred from the context.

The proof that Assumption 1 is satisfied is called the \textit{induction base}
(or \textit{base case}) of the proof. The proof that Assumption 2 is satisfied
is called the \textit{induction step} of the proof.

In order to prove that Assumption 2 is satisfied, you will usually want to fix
a $k\in\mathbb{Z}_{\geq g+1}$ such that $\mathcal{A}\left(  k-1\right)  $
holds, and then prove that $\mathcal{A}\left(  k\right)  $ holds. In other
words, you will usually want to fix $k\in\mathbb{Z}_{\geq g+1}$, assume that
$\mathcal{A}\left(  k-1\right)  $ holds, and then prove that $\mathcal{A}%
\left(  k\right)  $ holds. When doing so, it is common to refer to the
assumption that $\mathcal{A}\left(  k-1\right)  $ holds as the
\textit{induction hypothesis} (or \textit{induction assumption}).
\end{convention}

This language is exactly the same that was introduced in Convention
\ref{conv.ind.IPglang} for proofs by \textquotedblleft
standard\textquotedblright\ induction starting at $g$. The only difference
between proofs that use Theorem \ref{thm.ind.IPg} and proofs that use Theorem
\ref{thm.ind.IPg-1} is that the induction step in the former proofs assumes
$\mathcal{A}\left(  m\right)  $ and proves $\mathcal{A}\left(  m+1\right)  $,
whereas the induction step in the latter proofs assumes $\mathcal{A}\left(
k-1\right)  $ and proves $\mathcal{A}\left(  k\right)  $. (Of course, the
letters \textquotedblleft$m$\textquotedblright\ and \textquotedblleft%
$k$\textquotedblright\ are not set in stone; any otherwise unused letters can
be used in their stead. Thus, what distinguishes proofs that use Theorem
\ref{thm.ind.IPg} from proofs that use Theorem \ref{thm.ind.IPg-1} is not the
letter they use, but the \textquotedblleft$+1$\textquotedblright\ versus the
\textquotedblleft$-1$\textquotedblright.)

Let us repeat the above proof of Proposition \ref{prop.ind.alt-harm} (or, more
precisely, its non-computational part) using this language:

\begin{proof}
[Proof of Proposition \ref{prop.ind.alt-harm} (second version).]We must prove
(\ref{eq.prop.ind.alt-harm.claim}) for every $n\in\mathbb{N}$. In other words,
we must prove (\ref{eq.prop.ind.alt-harm.claim}) for every $n\in
\mathbb{Z}_{\geq0}$ (since $\mathbb{N}=\mathbb{Z}_{\geq0}$). We shall prove
this by induction on $n$ starting at $0$:

\textit{Induction base:} We have $\sum_{i=1}^{2\cdot0}\dfrac{\left(
-1\right)  ^{i-1}}{i}=\left(  \text{empty sum}\right)  =0$. Comparing this
with $\sum_{i=0+1}^{2\cdot0}\dfrac{1}{i}=\left(  \text{empty sum}\right)  =0$,
we obtain $\sum_{i=1}^{2\cdot0}\dfrac{\left(  -1\right)  ^{i-1}}{i}%
=\sum_{i=0+1}^{2\cdot0}\dfrac{1}{i}$. In other words,
(\ref{eq.prop.ind.alt-harm.claim}) holds for $n=0$. This completes the
induction base.

\textit{Induction step:} Let $k\in\mathbb{Z}_{\geq1}$. Assume that
(\ref{eq.prop.ind.alt-harm.claim}) holds for $n=k-1$. We must show that
(\ref{eq.prop.ind.alt-harm.claim}) holds for $n=k$.

We have $k\in\mathbb{Z}_{\geq1}$. In other words, $k$ is an integer and
satisfies $k\geq1$.

We have assumed that (\ref{eq.prop.ind.alt-harm.claim}) holds for $n=k-1$. In
other words,
\begin{equation}
\sum_{i=1}^{2\left(  k-1\right)  }\dfrac{\left(  -1\right)  ^{i-1}}{i}%
=\sum_{i=\left(  k-1\right)  +1}^{2\left(  k-1\right)  }\dfrac{1}{i}.
\label{pf.prop.ind.alt-harm.ver2.1}%
\end{equation}
From here, we can obtain%
\begin{equation}
\sum_{i=1}^{2k}\dfrac{\left(  -1\right)  ^{i-1}}{i}=\sum_{i=k+1}^{2k}\dfrac
{1}{i}. \label{pf.prop.ind.alt-harm.ver2.4}%
\end{equation}
(Indeed, we can derive (\ref{pf.prop.ind.alt-harm.ver2.4}) from
(\ref{pf.prop.ind.alt-harm.ver2.1}) in exactly the same way as we derived
(\ref{pf.prop.ind.alt-harm.8}) from (\ref{pf.prop.ind.alt-harm.IH}) in the
above first version of the proof of Proposition \ref{prop.ind.alt-harm};
nothing about this argument needs to be changed, so we have no reason to
repeat it.)

But the equality (\ref{pf.prop.ind.alt-harm.ver2.4}) shows that
(\ref{eq.prop.ind.alt-harm.claim}) holds for $n=k$. This completes the
induction step. Hence, (\ref{eq.prop.ind.alt-harm.claim}) is proven by
induction. This proves Proposition \ref{prop.ind.alt-harm}.
\end{proof}

\begin{exercise}
\label{exe.sum--2choosek-cases}Let $n\in\mathbb{N}$. Prove that%
\[
\sum_{k=0}^{n}\left(  -1\right)  ^{k}\left(  k+1\right)  =%
\begin{cases}
n/2+1, & \text{if }n\text{ is even};\\
-\left(  n+1\right)  /2, & \text{if }n\text{ is odd}%
\end{cases}
.
\]

\end{exercise}

The claim of Exercise \ref{exe.sum--2choosek-cases} can be rewritten as%
\[
1-2+3-4\pm\cdots+\left(  -1\right)  ^{n}\left(  n+1\right)  =%
\begin{cases}
n/2+1, & \text{if }n\text{ is even};\\
-\left(  n+1\right)  /2, & \text{if }n\text{ is odd}%
\end{cases}
.
\]

\section{\label{chp.binom}On binomial coefficients}

The present chapter is about \textit{binomial coefficients}. They are used in
almost every part of mathematics, and studying them provides good
opportunities to practice the arts of mathematical induction and of finding
combinatorial bijections.

Identities involving binomial coefficients are legion, and books have been
written about them (let me mention \cite[Chapter 5]{GKP} as a highly readable
introduction; but, e.g., \href{http://www.math.wvu.edu/~gould/}{Henry W.
Gould's website} goes far deeper down the rabbit hole). We shall only study a
few of these identities.

\subsection{\label{sect.binom.def}Definitions and basic properties}

\subsubsection{The definition}

Let us first define binomial coefficients:

\begin{definition}
\label{def.binom}Let $n\in\mathbb{N}$ and $m\in\mathbb{Q}$. Recall that $n!$
is a positive integer; thus, $n!\neq0$. (Keep in mind that $0!=1$.)

We define a rational number $\dbinom{m}{n}$ by%
\begin{equation}
\dbinom{m}{n}=\dfrac{m\left(  m-1\right)  \cdots\left(  m-n+1\right)  }{n!}.
\label{eq.binom.mn}%
\end{equation}
(This fraction is well-defined, since $n!\neq0$. When $n=0$, the numerator of
this fraction (i.e., the product $m\left(  m-1\right)  \cdots\left(
m-n+1\right)  $) is an empty product. Recall that, by convention, an empty
product is always defined to be $1$.)

This number $\dbinom{m}{n}$ is called a \textit{binomial coefficient}, and is
often pronounced \textquotedblleft$m$ choose $n$\textquotedblright.

We can extend this definition to the case when $m\in\mathbb{R}$ or
$m\in\mathbb{C}$ (rather than $m\in\mathbb{Q}$) by using the same equality
(\ref{eq.binom.mn}). Of course, in that case, $\dbinom{m}{n}$ will not be a
rational number anymore.
\end{definition}

\begin{example}
The formula (\ref{eq.binom.mn}) yields%
\begin{align*}
\dbinom{4}{2}  &  =\dfrac{4\left(  4-1\right)  }{2!}=\dfrac{4\left(
4-1\right)  }{2}=6;\\
\dbinom{5}{1}  &  =\dfrac{5}{1!}=\dfrac{5}{1}=5;\\
\dbinom{8}{3}  &  =\dfrac{8\left(  8-1\right)  \left(  8-2\right)  }%
{3!}=\dfrac{8\left(  8-1\right)  \left(  8-2\right)  }{6}=56.
\end{align*}

\end{example}

Here is a table of the binomial coefficients $\dbinom{n}{k}$ for all values
$n\in\left\{  -3,-2,-1,\ldots,6\right\}  $ and some of the values
$k\in\left\{  0,1,2,3,4,5\right\}  $. In the following table, each row
corresponds to a value of $n$, while each southwest-northeast diagonal
corresponds to a value of $k$:%
\[
\hspace{-0.3cm}%
\begin{tabular}
[c]{ll|ccccccccccccccccc}%
\vphantom{$g^b$} &  & $\phantom{20}$ & $\phantom{20}$ & $\phantom{20}$ &
$\phantom{20}$ &  &  &  &  &  &  & $\overset{k=0}{\swarrow}$ & $\phantom{20}$
& $\overset{k=1}{\swarrow}$ & $\phantom{20}$ & $\overset{k=2}{\swarrow}$ &
$\phantom{20}$ & $\overset{k=3}{\swarrow}$\\\hline
\vphantom{$g^b$} $n=-3$ & $\rightarrow$ &  &  &  &  &  &  &  &  &  & $1$ &  &
$-3$ &  & $6$ &  & $-10$ & \\
\vphantom{$g^b$} $n=-2$ & $\rightarrow$ &  &  &  &  &  &  &  &  & $1$ &  &
$-2$ &  & $3$ &  & $-4$ &  & \\
\vphantom{$g^b$} $n=-1$ & $\rightarrow$ &  &  &  &  &  &  &  & $1$ &  & $-1$ &
& $1$ &  & $-1$ &  & $1$ & \\
\vphantom{$g^b$} $n=0$ & $\rightarrow$ &  &  &  &  &  &  & $1$ &  & $0$ &  &
$0$ &  & $0$ &  & $0$ &  & \\
\vphantom{$g^b$} $n=1$ & $\rightarrow$ &  &  &  &  &  & $1$ &  & $1$ &  & $0$
&  & $0$ &  & $0$ &  & $0$ & \\
\vphantom{$g^b$} $n=2$ & $\rightarrow$ &  &  &  &  & $1$ &  & $2$ &  & $1$ &
& $0$ &  & $0$ &  & $0$ &  & \\
\vphantom{$g^b$} $n=3$ & $\rightarrow$ &  &  &  & $1$ &  & $3$ &  & $3$ &  &
$1$ &  & $0$ &  & $0$ &  & $0$ & \\
\vphantom{$g^b$} $n=4$ & $\rightarrow$ &  &  & $1$ &  & $4$ &  & $6$ &  & $4$
&  & $1$ &  & $0$ &  & $0$ &  & \\
\vphantom{$g^b$} $n=5$ & $\rightarrow$ &  & $1$ &  & $5$ &  & $10$ &  & $10$ &
& $5$ &  & $1$ &  & $0$ &  & $0$ & \\
\vphantom{$g^b$} $n=6$ & $\rightarrow$ & $1$ &  & $6$ &  & $15$ &  & $20$ &  &
$15$ &  & $6$ &  & $1$ &  & $0$ &  &
\end{tabular}
\ \ \ \ \ \ \
\]

The binomial coefficients $\dbinom{m}{n}$ form the so-called
\textit{\href{https://en.wikipedia.org/wiki/Pascal's_triangle}{\textit{Pascal's
triangle}}}\footnote{More precisely, the numbers $\dbinom{m}{n}$ for
$m\in\mathbb{N}$ and $n\in\left\{  0,1,\ldots,m\right\}  $ form Pascal's
triangle. Nevertheless, the \textquotedblleft other\textquotedblright%
\ binomial coefficients (particularly the ones where $m$ is a negative
integer) are highly useful, too.}. Let us state a few basic properties of
these numbers:

\subsubsection{Simple formulas}

\begin{proposition}
\label{prop.binom.00}Let $m\in\mathbb{Q}$.

\textbf{(a)} We have%
\begin{equation}
\dbinom{m}{0}=1. \label{eq.binom.00}%
\end{equation}

\textbf{(b)} We have%
\begin{equation}
\dbinom{m}{1}=m. \label{eq.binom.01}%
\end{equation}

\end{proposition}

\begin{proof}
[Proof of Proposition \ref{prop.binom.00}.]\textbf{(a)} The definition of
$\dbinom{m}{0}$ yields
\[
\dbinom{m}{0}=\dfrac{m\left(  m-1\right)  \cdots\left(  m-0+1\right)  }{0!}.
\]
Since $m\left(  m-1\right)  \cdots\left(  m-0+1\right)  =\left(  \text{a
product of }0\text{ integers}\right)  =1$, this rewrites as $\dbinom{m}%
{0}=\dfrac{1}{0!}=1$ (since $0!=1$). This proves Proposition
\ref{prop.binom.00} \textbf{(a)}.

\textbf{(b)} The definition of $\dbinom{m}{1}$ yields%
\[
\dbinom{m}{1}=\dfrac{m\left(  m-1\right)  \cdots\left(  m-1+1\right)  }{1!}.
\]
Since $m\left(  m-1\right)  \cdots\left(  m-1+1\right)  =m$ (because the
product $m\left(  m-1\right)  \cdots\left(  m-1+1\right)  $ consists of $1$
factor only, and this factor is $m$), this rewrites as $\dbinom{m}{1}%
=\dfrac{m}{1!}=m$ (since $1!=1$). This proves Proposition \ref{prop.binom.00}
\textbf{(b)}.
\end{proof}

\begin{proposition}
\label{prop.binom.formula}Let $m\in\mathbb{N}$ and $n\in\mathbb{N}$ be such
that $m\geq n$. Then,%
\begin{equation}
\dbinom{m}{n}=\dfrac{m!}{n!\left(  m-n\right)  !}. \label{eq.binom.formula}%
\end{equation}

\end{proposition}

\begin{remark}
\textbf{Caution:} The formula (\ref{eq.binom.formula}) holds only for
$m\in\mathbb{N}$ and $n\in\mathbb{N}$ satisfying $m\geq n$. Thus, neither
$\dbinom{-3}{2}$ nor $\dbinom{1/3}{3}$ nor $\dbinom{2}{5}$ can be computed
using this formula! Definition \ref{def.binom} thus can be used to compute
$\dbinom{m}{n}$ in many more cases than (\ref{eq.binom.formula}) does.
\end{remark}

\begin{proof}
[Proof of Proposition \ref{prop.binom.formula}.]Multiplying both sides of the
equality (\ref{eq.binom.mn}) with $n!$, we obtain
\begin{equation}
n!\cdot\dbinom{m}{n}=m\left(  m-1\right)  \cdots\left(  m-n+1\right)  .
\label{pf.prop.binom.formula.1}%
\end{equation}
But%
\begin{align*}
m!  &  =m\left(  m-1\right)  \cdots1=\left(  m\left(  m-1\right)
\cdots\left(  m-n+1\right)  \right)  \cdot\underbrace{\left(  \left(
m-n\right)  \left(  m-n-1\right)  \cdots1\right)  }_{=\left(  m-n\right)  !}\\
&  =\left(  m\left(  m-1\right)  \cdots\left(  m-n+1\right)  \right)
\cdot\left(  m-n\right)  !,
\end{align*}
so that $\dfrac{m!}{\left(  m-n\right)  !}=m\left(  m-1\right)  \cdots\left(
m-n+1\right)  $. Comparing this with (\ref{pf.prop.binom.formula.1}), we
obtain $n!\cdot\dbinom{m}{n}=\dfrac{m!}{\left(  m-n\right)  !}$. Dividing this
equality by $n!$, we obtain $\dbinom{m}{n}=\dfrac{m!}{n!\left(  m-n\right)
!}$. Thus, Proposition \ref{prop.binom.formula} is proven.
\end{proof}

\begin{proposition}
\label{prop.binom.0}Let $m\in\mathbb{N}$ and $n\in\mathbb{N}$ be such that
$m<n$. Then,
\begin{equation}
\dbinom{m}{n}=0. \label{eq.binom.0}%
\end{equation}

\end{proposition}

\begin{remark}
\textbf{Caution:} The formula (\ref{eq.binom.0}) is not true if we drop the
condition $m\in\mathbb{N}$. For example, $\dbinom{-3}{2}=6\neq0$ despite
$-3<2$.
\end{remark}

\begin{proof}
[Proof of Proposition \ref{prop.binom.0}.]We have $m\geq0$ (since
$m\in\mathbb{N}$). Also, $m<n$, so that $m\leq n-1$ (since $m$ and $n$ are
integers). Thus, $m\in\left\{  0,1,\ldots,n-1\right\}  $. Hence, $m-m$ is one
of the $n$ integers $m,m-1,\ldots,m-n+1$. Thus, one of the $n$ factors of the
product $m\left(  m-1\right)  \cdots\left(  m-n+1\right)  $ is $m-m=0$.
Therefore, the whole product $m\left(  m-1\right)  \cdots\left(  m-n+1\right)
$ is $0$ (because if one of the factors of a product is $0$, then the whole
product must be $0$). Thus, $m\left(  m-1\right)  \cdots\left(  m-n+1\right)
=0$. Hence, (\ref{eq.binom.mn}) becomes%
\begin{align*}
\dbinom{m}{n}  &  =\dfrac{m\left(  m-1\right)  \cdots\left(  m-n+1\right)
}{n!}=\dfrac{0}{n!}\ \ \ \ \ \ \ \ \ \ \left(  \text{since }m\left(
m-1\right)  \cdots\left(  m-n+1\right)  =0\right) \\
&  =0.
\end{align*}
This proves Proposition \ref{prop.binom.0}.
\end{proof}

\begin{proposition}
\label{prop.binom.symm}Let $m\in\mathbb{N}$ and $n\in\mathbb{N}$ be such that
$m\geq n$. Then,%
\begin{equation}
\dbinom{m}{n}=\dbinom{m}{m-n}. \label{eq.binom.symm}%
\end{equation}

\end{proposition}

Proposition \ref{prop.binom.symm} is commonly known as the \textit{symmetry
identity for the binomial coefficients}. Notice that Proposition
\ref{prop.binom.symm} becomes false (and, with our definitions, actually
meaningless) if the requirement that $m\in\mathbb{N}$ is dropped.

\begin{proof}
[Proof of Proposition \ref{prop.binom.symm}.]We have $m-n\in\mathbb{N}$ (since
$m\geq n$) and $m\geq m-n$ (since $n\geq0$ (since $n\in\mathbb{N}$)). Hence,
(\ref{eq.binom.formula}) (applied to $m-n$ instead of $n$) yields
\[
\dbinom{m}{m-n}=\dfrac{m!}{\left(  m-n\right)  !\left(  m-\left(  m-n\right)
\right)  !}=\dfrac{m!}{\left(  m-\left(  m-n\right)  \right)  !\left(
m-n\right)  !}=\dfrac{m!}{n!\left(  m-n\right)  !}%
\]
(since $m-\left(  m-n\right)  =n$). Compared with (\ref{eq.binom.formula}),
this yields $\dbinom{m}{n}=\dbinom{m}{m-n}$. Proposition \ref{prop.binom.symm}
is thus proven.
\end{proof}

\begin{proposition}
\label{prop.binom.mm}Let $m\in\mathbb{N}$. Then,%
\begin{equation}
\dbinom{m}{m}=1. \label{eq.binom.mm}%
\end{equation}

\end{proposition}

\begin{proof}
[Proof of Proposition \ref{prop.binom.mm}.]The equality (\ref{eq.binom.symm})
(applied to $n=m$) yields $\dbinom{m}{m}=\dbinom{m}{m-m}=\dbinom{m}{0}=1$
(according to (\ref{eq.binom.00})). This proves Proposition
\ref{prop.binom.mm}.
\end{proof}

\begin{exercise}
\label{exe.multinom1}Let $m\in\mathbb{N}$ and $\left(  k_{1},k_{2}%
,\ldots,k_{m}\right)  \in\mathbb{N}^{m}$. Prove that $\dfrac{\left(
k_{1}+k_{2}+\cdots+k_{m}\right)  !}{k_{1}!k_{2}!\cdots k_{m}!}$ is a positive integer.
\end{exercise}

\begin{remark}
\label{rmk.multinom1}Let $m\in\mathbb{N}$ and $\left(  k_{1},k_{2}%
,\ldots,k_{m}\right)  \in\mathbb{N}^{m}$. Exercise \ref{exe.multinom1} shows
that $\dfrac{\left(  k_{1}+k_{2}+\cdots+k_{m}\right)  !}{k_{1}!k_{2}!\cdots
k_{m}!}$ is a positive integer. This positive integer is called a
\textit{multinomial coefficient}, and is often denoted by $\dbinom{n}%
{k_{1},k_{2},\ldots,k_{m}}$, where $n=k_{1}+k_{2}+\cdots+k_{m}$. (We shall
avoid this particular notation, since it makes the meaning of $\dbinom{n}{k}$
slightly ambiguous: It could mean both the binomial coefficient $\dbinom{n}%
{k}$ and the multinomial coefficient $\dbinom{n}{k_{1},k_{2},\ldots,k_{m}}$
for $\left(  k_{1},k_{2},\ldots,k_{m}\right)  =\left(  k\right)  $.
Fortunately, the ambiguity is not really an issue, because the only situation
in which both meanings make sense is when $k=n\in\mathbb{N}$, but in this case
both interpretations give the same value $1$.)
\end{remark}

\begin{exercise}
\label{exe.bin.-1/2}Let $n\in\mathbb{N}$.

\textbf{(a)} Prove that%
\[
\left(  2n-1\right)  \cdot\left(  2n-3\right)  \cdot\cdots\cdot1=\dfrac
{\left(  2n\right)  !}{2^{n}n!}.
\]
(The left hand side is understood to be the product of all odd integers from
$1$ to $2n-1$.)

\textbf{(b)} Prove that
\[
\dbinom{-1/2}{n}=\left(  \dfrac{-1}{4}\right)  ^{n}\dbinom{2n}{n}.
\]

\textbf{(c)} Prove that%
\[
\dbinom{-1/3}{n}\dbinom{-2/3}{n}=\dfrac{\left(  3n\right)  !}{\left(
3^{n}n!\right)  ^{3}}.
\]

\end{exercise}

\subsubsection{The recurrence relation of the binomial coefficients}

\begin{proposition}
\label{prop.binom.rec}Let $m\in\mathbb{Q}$ and $n\in\left\{  1,2,3,\ldots
\right\}  $. Then,%
\begin{equation}
\dbinom{m}{n}=\dbinom{m-1}{n-1}+\dbinom{m-1}{n}. \label{eq.binom.rec.m}%
\end{equation}

\end{proposition}

\begin{proof}
[Proof of Proposition \ref{prop.binom.rec}.]From $n\in\left\{  1,2,3,\ldots
\right\}  $, we obtain $n!=n\cdot\left(  n-1\right)  !$, so that $\left(
n-1\right)  !=n!/n$ and thus $\dfrac{1}{\left(  n-1\right)  !}=\dfrac{1}%
{n!/n}=\dfrac{1}{n!}\cdot n$.

The definition of $\dbinom{m}{n-1}$ yields
\begin{align*}
\dbinom{m}{n-1}  &  =\dfrac{m\left(  m-1\right)  \cdots\left(  m-\left(
n-1\right)  +1\right)  }{\left(  n-1\right)  !}\\
&  =\dfrac{1}{\left(  n-1\right)  !}\cdot\left(  m\left(  m-1\right)
\cdots\left(  m-\left(  n-1\right)  +1\right)  \right)  .
\end{align*}
The same argument (applied to $m-1$ instead of $m$) yields%
\begin{align}
\dbinom{m-1}{n-1}  &  =\underbrace{\dfrac{1}{\left(  n-1\right)  !}}%
_{=\dfrac{1}{n!}\cdot n}\cdot\left(  \left(  m-1\right)  \underbrace{\left(
\left(  m-1\right)  -1\right)  }_{=m-2}\cdots\underbrace{\left(  \left(
m-1\right)  -\left(  n-1\right)  +1\right)  }_{=m-n+1}\right) \nonumber\\
&  =\dfrac{1}{n!}\cdot n\cdot\left(  \left(  m-1\right)  \left(  m-2\right)
\cdots\left(  m-n+1\right)  \right)  . \label{eq.binom.rec.pf.1}%
\end{align}

On the other hand,
\[
\dbinom{m}{n}=\dfrac{m\left(  m-1\right)  \cdots\left(  m-n+1\right)  }%
{n!}=\dfrac{1}{n!}\left(  m\left(  m-1\right)  \cdots\left(  m-n+1\right)
\right)  .
\]
The same argument (applied to $m-1$ instead of $m$) yields%
\begin{align*}
\dbinom{m-1}{n}  &  =\dfrac{1}{n!}\left(  \left(  m-1\right)
\underbrace{\left(  \left(  m-1\right)  -1\right)  }_{=m-2}\cdots
\underbrace{\left(  \left(  m-1\right)  -n+1\right)  }_{=m-n}\right) \\
&  =\dfrac{1}{n!}\underbrace{\left(  \left(  m-1\right)  \left(  m-2\right)
\cdots\left(  m-n\right)  \right)  }_{=\left(  \left(  m-1\right)  \left(
m-2\right)  \cdots\left(  m-n+1\right)  \right)  \cdot\left(  m-n\right)  }\\
&  =\dfrac{1}{n!}\left(  \left(  m-1\right)  \left(  m-2\right)  \cdots\left(
m-n+1\right)  \right)  \cdot\left(  m-n\right) \\
&  =\dfrac{1}{n!}\left(  m-n\right)  \cdot\left(  \left(  m-1\right)  \left(
m-2\right)  \cdots\left(  m-n+1\right)  \right)  .
\end{align*}
Adding (\ref{eq.binom.rec.pf.1}) to this equality, we obtain%
\begin{align*}
&  \dbinom{m-1}{n}+\dbinom{m-1}{n-1}\\
&  =\dfrac{1}{n!}\left(  m-n\right)  \cdot\left(  \left(  m-1\right)  \left(
m-2\right)  \cdots\left(  m-n+1\right)  \right) \\
&  \ \ \ \ \ \ \ \ \ \ +\dfrac{1}{n!}\cdot n\cdot\left(  \left(  m-1\right)
\left(  m-2\right)  \cdots\left(  m-n+1\right)  \right) \\
&  =\dfrac{1}{n!}\underbrace{\left(  \left(  m-n\right)  +n\right)  }%
_{=m}\cdot\left(  \left(  m-1\right)  \left(  m-2\right)  \cdots\left(
m-n+1\right)  \right) \\
&  =\dfrac{1}{n!}\underbrace{m\cdot\left(  \left(  m-1\right)  \left(
m-2\right)  \cdots\left(  m-n+1\right)  \right)  }_{=m\left(  m-1\right)
\cdots\left(  m-n+1\right)  }=\dfrac{1}{n!}\left(  m\left(  m-1\right)
\cdots\left(  m-n+1\right)  \right) \\
&  =\dbinom{m}{n}\ \ \ \ \ \ \ \ \ \ \left(  \text{since }\dbinom{m}{n}%
=\dfrac{1}{n!}\left(  m\left(  m-1\right)  \cdots\left(  m-n+1\right)
\right)  \right)  .
\end{align*}
This proves Proposition \ref{prop.binom.rec}.
\end{proof}

The formula (\ref{eq.binom.rec.m}) is known as the \textit{recurrence relation
of the binomial coefficients}\footnote{Often it is extended to the case $n=0$
by setting $\dbinom{m}{-1}=0$. It then follows from (\ref{eq.binom.00}) in
this case.
\par
The formula (\ref{eq.binom.rec.m}) is responsible for the fact that
\textquotedblleft every number in Pascal's triangle is the sum of the two
numbers above it\textquotedblright. (Of course, if you use this fact as a
\textit{definition} of Pascal's triangle, then (\ref{eq.binom.rec.m}) is
conversely responsible for the fact that the numbers in this triangle are the
binomial coefficients.)}.

\begin{exercise}
\label{exe.binom.hockey1}\textbf{(a)} Prove that every $n\in\mathbb{N}$ and
$q\in\mathbb{Q}$ satisfy%
\[
\sum_{r=0}^{n}\dbinom{r+q}{r}=\dbinom{n+q+1}{n}.
\]

\textbf{(b)} Prove that every $n\in\left\{  -1,0,1,\ldots\right\}  $ and
$k\in\mathbb{N}$ satisfy%
\[
\sum_{i=0}^{n}\dbinom{i}{k}=\sum_{i=k}^{n}\dbinom{i}{k}=\dbinom{n+1}{k+1}.
\]
(Keep in mind that $\sum_{i=k}^{n}\dbinom{i}{k}$ is an empty sum whenever
$n<k$.)
\end{exercise}

The claim of Exercise \ref{exe.binom.hockey1} \textbf{(b)} is one of several
formulas known as the
\textit{\href{https://en.wikipedia.org/wiki/Hockey-stick_identity}{\textit{hockey-stick
identity}}} (due to the fact that marking the binomial coefficients appearing
in it in Pascal's triangle results in a shape resembling a hockey
stick\footnote{See \url{https://math.stackexchange.com/q/1490794} for an
illustration.}); it appears, e.g., in \cite[Identity 11.10]{Galvin} (or,
rather, the second equality sign of Exercise \ref{exe.binom.hockey1}
\textbf{(b)} appears there, but the rest is easy).

\subsubsection{The combinatorial interpretation of binomial coefficients}

\begin{proposition}
\label{prop.binom.subsets}If $m\in\mathbb{N}$ and $n\in\mathbb{N}$, and if $S$
is an $m$-element set, then
\begin{equation}
\dbinom{m}{n}\text{ is the number of all }n\text{-element subsets of
}S\text{.} \label{eq.binom.subsets}%
\end{equation}

\end{proposition}

In less formal terms, Proposition \ref{prop.binom.subsets} says the following:
If $m\in\mathbb{N}$ and $n\in\mathbb{N}$, then $\dbinom{m}{n}$ is the number
of ways to pick out $n$ among $m$ given objects, without
replacement\footnote{That is, one must not pick out the same object twice.}
and without regard for the order in which they are picked out. (Probabilists
call this \textquotedblleft unordered samples without
replacement\textquotedblright.)

\begin{example}
Proposition \ref{prop.binom.subsets} (applied to $m=4$, $n=2$ and $S=\left\{
0,1,2,3\right\}  $) shows that $\dbinom{4}{2}$ is the number of all
$2$-element subsets of $\left\{  0,1,2,3\right\}  $ (since $\left\{
0,1,2,3\right\}  $ is a $4$-element set). And indeed, this is easy to verify
by brute force: The $2$-element subsets of $\left\{  0,1,2,3\right\}  $ are%
\[
\left\{  0,1\right\}  ,\ \left\{  0,2\right\}  ,\ \left\{  0,3\right\}
,\ \left\{  1,2\right\}  ,\ \left\{  1,3\right\}  \text{ and }\left\{
2,3\right\}  ,
\]
so there are $6=\dbinom{4}{2}$ of them.
\end{example}

\begin{remark}
\textbf{Caution:} Proposition \ref{prop.binom.subsets} says nothing about
binomial coefficients $\dbinom{m}{n}$ with negative $m$. Indeed, there are no
$m$-element sets $S$ when $m$ is negative; thus, Proposition
\ref{prop.binom.subsets} would be vacuously true\footnotemark\ when $m$ is
negative, but this would not help us computing binomial coefficients
$\dbinom{m}{n}$ with negative $m$.

Actually, when $m\in\mathbb{Z}$ is negative, the number $\dbinom{m}{n}$ is
positive for $n$ even and negative for $n$ odd (easy exercise), and so an
interpretation of $\dbinom{m}{n}$ as a number of ways to do something is
rather unlikely. (On the other hand, $\left(  -1\right)  ^{n}\dbinom{m}{n}$
does have such an interpretation.)
\end{remark}

\footnotetext{Recall that a mathematical statement of the form
\textquotedblleft if $\mathcal{A}$, then $\mathcal{B}$\textquotedblright\ is
said to be \textit{vacuously true} if $\mathcal{A}$ never holds. For example,
the statement \textquotedblleft if $0=1$, then every integer is
odd\textquotedblright\ is vacuously true, because $0=1$ is false. Proposition
\ref{prop.binom.subsets} is vacuously true when $m$ is negative, because the
condition \textquotedblleft$S$ is an $m$-element set\textquotedblright\ never
holds when $m$ is negative.
\par
By the laws of logic, a vacuously true statement is always true! See
Convention \ref{conv.logic.vacuous} for a discussion of this principle.}

\begin{remark}
Some authors (for example, those of \cite{LeLeMe16} and of \cite{Galvin}) use
(\ref{eq.binom.subsets}) as the \textit{definition} of $\dbinom{m}{n}$. This
is a legitimate definition of $\dbinom{m}{n}$ in the case when $m$ and $n$ are
nonnegative integers (and, of course, equivalent to our definition); but it is
not as general as ours, since it does not extend to negative (or non-integer)
$m$.
\end{remark}

\begin{exercise}
\label{exe.prop.binom.subsets}Prove Proposition \ref{prop.binom.subsets}.
\end{exercise}

Proposition \ref{prop.binom.subsets} is one of the most basic facts of
\textit{enumerative combinatorics} -- the part of mathematics that is mostly
concerned with counting problems (i.e., the study of the sizes of finite
sets). We will encounter some further results from enumerative combinatorics
below (e.g., Exercise \ref{exe.multichoose} and Exercise \ref{exe.ps2.2.3});
but we shall not go deep into this subject. More serious expositions of
enumerative combinatorics include Loehr's textbook \cite{Loehr-BC}, Galvin's
lecture notes \cite{Galvin}, Aigner's book \cite{Aigner07}, and Stanley's
two-volume treatise (\cite{Stanley-EC1} and \cite{Stanley-EC2}).

\subsubsection{Upper negation}

\begin{proposition}
\label{prop.binom.upper-neg}Let $m\in\mathbb{Q}$ and $n\in\mathbb{N}$. Then,%
\begin{equation}
\dbinom{m}{n}=\left(  -1\right)  ^{n}\dbinom{n-m-1}{n}.
\label{eq.binom.upper-neg}%
\end{equation}

\end{proposition}

\begin{proof}
[Proof of Proposition \ref{prop.binom.upper-neg}.]The equality
(\ref{eq.binom.mn}) (applied to $n-m-1$ instead of $m$) yields%
\begin{align*}
\dbinom{n-m-1}{n}  &  =\dfrac{\left(  n-m-1\right)  \left(  \left(
n-m-1\right)  -1\right)  \cdots\left(  \left(  n-m-1\right)  -n+1\right)
}{n!}\\
&  =\dfrac{1}{n!}\left(  n-m-1\right)  \left(  \left(  n-m-1\right)
-1\right)  \cdots\underbrace{\left(  \left(  n-m-1\right)  -n+1\right)
}_{=-m}\\
&  =\dfrac{1}{n!}\underbrace{\left(  n-m-1\right)  \left(  \left(
n-m-1\right)  -1\right)  \cdots\left(  -m\right)  }_{\substack{=\left(
-m\right)  \left(  -m+1\right)  \cdots\left(  n-m-1\right)  \\\text{(here, we
have just reversed the order of the factors in the product)}}}\\
&  =\dfrac{1}{n!}\underbrace{\left(  -m\right)  }_{=\left(  -1\right)
m}\underbrace{\left(  -m+1\right)  }_{=\left(  -1\right)  \left(  m-1\right)
}\cdots\underbrace{\left(  n-m-1\right)  }_{=\left(  -1\right)  \left(
m-n+1\right)  }\\
&  =\dfrac{1}{n!}\left(  \left(  -1\right)  m\right)  \left(  \left(
-1\right)  \left(  m-1\right)  \right)  \cdots\left(  \left(  -1\right)
\left(  m-n+1\right)  \right) \\
&  =\dfrac{1}{n!}\left(  -1\right)  ^{n}\left(  m\left(  m-1\right)
\cdots\left(  m-n+1\right)  \right)  ,
\end{align*}
so that%
\begin{align*}
\left(  -1\right)  ^{n}\dbinom{n-m-1}{n}  &  =\left(  -1\right)  ^{n}%
\cdot\dfrac{1}{n!}\left(  -1\right)  ^{n}\left(  m\left(  m-1\right)
\cdots\left(  m-n+1\right)  \right) \\
&  =\underbrace{\left(  -1\right)  ^{n}\left(  -1\right)  ^{n}}%
_{\substack{=\left(  -1\right)  ^{n+n}=\left(  -1\right)  ^{2n}%
=1\\\text{(since }2n\text{ is even)}}}\cdot\dfrac{1}{n!}\left(  m\left(
m-1\right)  \cdots\left(  m-n+1\right)  \right) \\
&  =\dfrac{1}{n!}\left(  m\left(  m-1\right)  \cdots\left(  m-n+1\right)
\right)  =\dfrac{m\left(  m-1\right)  \cdots\left(  m-n+1\right)  }{n!}.
\end{align*}
Compared with (\ref{eq.binom.mn}), this yields $\dbinom{m}{n}=\left(
-1\right)  ^{n}\dbinom{n-m-1}{n}$. Proposition \ref{prop.binom.upper-neg} is
therefore proven.
\end{proof}

The formula (\ref{eq.binom.upper-neg}) is known as the \textit{upper negation
formula}.

\begin{corollary}
\label{cor.binom.-1}Let $n\in\mathbb{N}$. Then,%
\[
\dbinom{-1}{n}=\left(  -1\right)  ^{n}.
\]

\end{corollary}

\begin{proof}
[Proof of Corollary \ref{cor.binom.-1}.]Proposition \ref{prop.binom.upper-neg}
(applied to $m=-1$) yields%
\begin{align*}
\dbinom{-1}{n}  &  =\left(  -1\right)  ^{n}\dbinom{n-\left(  -1\right)  -1}%
{n}=\left(  -1\right)  ^{n}\underbrace{\dbinom{n}{n}}_{\substack{=1\\\text{(by
Proposition \ref{prop.binom.mm}}\\\text{(applied to }m=n\text{))}}}\\
&  \ \ \ \ \ \ \ \ \ \ \left(  \text{since }n-\left(  -1\right)  -1=n\right)
\\
&  =\left(  -1\right)  ^{n}.
\end{align*}
This proves Corollary \ref{cor.binom.-1}.
\end{proof}

\begin{exercise}
\label{exe.bin.-1and-2}\textbf{(a)} Show that $\dbinom{-1}{k}=\left(
-1\right)  ^{k}$ for each $k\in\mathbb{N}$.

\textbf{(b)} Show that $\dbinom{-2}{k}=\left(  -1\right)  ^{k}\left(
k+1\right)  $ for each $k\in\mathbb{N}$.

\textbf{(c)} Show that $\dfrac{1!\cdot2!\cdot\cdots\cdot\left(  2n\right)
!}{n!}=2^{n}\cdot\left(  \prod_{i=1}^{n}\left(  \left(  2i-1\right)  !\right)
\right)  ^{2}$ for each $n\in\mathbb{N}$.
\end{exercise}

\begin{remark}
Parts \textbf{(a)} and \textbf{(b)} of Exercise \ref{exe.bin.-1and-2} are
known facts (actually, part \textbf{(a)} is just a repetition of Corollary
\ref{cor.binom.-1}, for the purpose of making the analogy to part \textbf{(b)}
more visible). Part \textbf{(c)} is a generalization of a puzzle posted on
\url{https://www.reddit.com/r/math/comments/7rybhp/factorial_problem/} . (The
puzzle boils down to the fact that $\dfrac{1!\cdot2!\cdot\cdots\cdot\left(
2n\right)  !}{n!}$ is a perfect square when $n\in\mathbb{N}$ is even. But this
follows from Exercise \ref{exe.bin.-1and-2} \textbf{(c)}, because when
$n\in\mathbb{N}$ is even, both factors $2^{n}$ and $\left(  \prod_{i=1}%
^{n}\left(  \left(  2i-1\right)  !\right)  \right)  ^{2}$ on the right hand
side of Exercise \ref{exe.bin.-1and-2} \textbf{(c)} are perfect squares.)
\end{remark}

\subsubsection{Binomial coefficients of integers are integers}

\begin{lemma}
\label{lem.binom.intN}Let $m\in\mathbb{N}$ and $n\in\mathbb{N}$. Then,
\[
\dbinom{m}{n}\in\mathbb{N}.
\]

\end{lemma}

\begin{proof}
[Proof of Lemma \ref{lem.binom.intN}.]We have $m\in\mathbb{N}$. Thus, there
exists an $m$-element set $S$ (for example, $S=\left\{  1,2,\ldots,m\right\}
$). Consider such an $S$. Then, $\dbinom{m}{n}$ is the number of all
$n$-element subsets of $S$ (because of (\ref{eq.binom.subsets})). Hence,
$\dbinom{m}{n}$ is a nonnegative integer, so that $\dbinom{m}{n}\in\mathbb{N}%
$. This proves Lemma \ref{lem.binom.intN}.
\end{proof}

It is also easy to prove Lemma \ref{lem.binom.intN} by induction on $m$, using
\eqref{eq.binom.00} and \eqref{eq.binom.0} in the induction base and using
(\ref{eq.binom.rec.m}) in the induction step.

\begin{proposition}
\label{prop.binom.int}Let $m\in\mathbb{Z}$ and $n\in\mathbb{N}$. Then,%
\begin{equation}
\dbinom{m}{n}\in\mathbb{Z}. \label{eq.binom.int}%
\end{equation}

\end{proposition}

\begin{proof}
[Proof of Proposition \ref{prop.binom.int}.]We need to show
(\ref{eq.binom.int}). We are in one of the following two cases:

\textit{Case 1:} We have $m\geq0$.

\textit{Case 2:} We have $m<0$.

Let us first consider Case 1. In this case, we have $m\geq0$. Hence,
$m\in\mathbb{N}$. Thus, Lemma \ref{lem.binom.intN} yields $\dbinom{m}{n}%
\in\mathbb{N}\subseteq\mathbb{Z}$. This proves (\ref{eq.binom.int}) in Case 1.

Let us now consider Case 2. In this case, we have $m<0$. Thus, $m\leq-1$
(since $m$ is an integer), so that $m+1\leq0$, so that
$n-m-1=n-\underbrace{\left(  m+1\right)  }_{\leq0}\geq n\geq0$. Hence,
$n-m-1\in\mathbb{N}$. Therefore, Lemma \ref{lem.binom.intN} (applied to
$n-m-1$ instead of $m$) yields $\dbinom{n-m-1}{n}\in\mathbb{N}\subseteq
\mathbb{Z}$. Now, (\ref{eq.binom.upper-neg}) shows that $\dbinom{m}%
{n}=\underbrace{\left(  -1\right)  ^{n}}_{\in\mathbb{Z}}\underbrace{\dbinom
{n-m-1}{n}}_{\in\mathbb{Z}}\in\mathbb{Z}$ (here, we have used the fact that
the product of two integers is an integer). This proves (\ref{eq.binom.int})
in Case 2.

We thus have proven (\ref{eq.binom.int}) in each of the two Cases 1 and 2. We
can therefore conclude that (\ref{eq.binom.int}) always holds. Thus,
Proposition \ref{prop.binom.int} is proven.
\end{proof}

The above proof of Proposition \ref{prop.binom.int} may well be the simplest
one. There is another proof, which uses Theorem \ref{thm.ind.IPg+-}, but it is
more complicated\footnote{It requires an induction on $n$ nested inside the
induction step of the induction on $m$.}. There is yet another proof using
basic number theory (specifically, checking how often a prime $p$ appears in
the numerator and the denominator of $\dbinom{m}{n}=\dfrac{m\left(
m-1\right)  \cdots\left(  m-n+1\right)  }{n!}$), but this is not quite easy.

\subsubsection{The binomial formula}

\begin{proposition}
\label{prop.binom.binomial}Let $x$ and $y$ be two rational numbers (or real
numbers, or complex numbers). Then,%
\begin{equation}
\left(  x+y\right)  ^{n}=\sum_{k=0}^{n}\dbinom{n}{k}x^{k}y^{n-k}
\label{eq.binom.binomial}%
\end{equation}
for every $n\in\mathbb{N}$.
\end{proposition}

Proposition \ref{prop.binom.binomial} is the famous \textit{binomial formula}
(also known as the \textit{binomial theorem}) and has a well-known standard
proof by induction over $n$ (using (\ref{eq.binom.rec.m}) and
(\ref{eq.binom.00}))\footnote{See Exercise \ref{exe.prop.binom.binomial} for
this proof.}. Some versions of it hold for negative $n$ as well (but not in
the exact form (\ref{eq.binom.binomial}), and not without restrictions).

\begin{exercise}
\label{exe.prop.binom.binomial}Prove Proposition \ref{prop.binom.binomial}.
\end{exercise}

There is an analogue of Proposition \ref{prop.binom.binomial} for a sum of $m$
rational numbers (rather than $2$ rational numbers); it is called the
\textquotedblleft multinomial formula\textquotedblright\ (and involves the
multinomial coefficients from Remark \ref{rmk.multinom1}). We shall state it
in a more general setting in Exercise \ref{exe.multinom2}.

\subsubsection{The absorption identity}

\begin{proposition}
\label{prop.binom.X-1}Let $n\in\left\{  1,2,3,\ldots\right\}  $ and
$m\in\mathbb{Q}$. Then,%
\begin{equation}
\dbinom{m}{n}=\dfrac{m}{n}\dbinom{m-1}{n-1}. \label{eq.binom.X-1}%
\end{equation}

\end{proposition}

\begin{proof}
[Proof of Proposition \ref{prop.binom.X-1}.]The definition of $\dbinom
{m-1}{n-1}$ yields%
\begin{align*}
\dbinom{m-1}{n-1}  &  =\dfrac{\left(  m-1\right)  \left(  \left(  m-1\right)
-1\right)  \cdots\left(  \left(  m-1\right)  -\left(  n-1\right)  +1\right)
}{\left(  n-1\right)  !}\\
&  =\dfrac{\left(  m-1\right)  \left(  m-2\right)  \cdots\left(  m-n+1\right)
}{\left(  n-1\right)  !}%
\end{align*}
(since $\left(  m-1\right)  -1=m-2$ and $\left(  m-1\right)  -\left(
n-1\right)  +1=m-n+1$). Multiplying both sides of this equality by $\dfrac
{m}{n}$, we obtain%
\begin{align*}
\dfrac{m}{n}\dbinom{m-1}{n-1}  &  =\dfrac{m}{n}\cdot\dfrac{\left(  m-1\right)
\left(  m-2\right)  \cdots\left(  m-n+1\right)  }{\left(  n-1\right)  !}\\
&  =\dfrac{m\left(  m-1\right)  \left(  m-2\right)  \cdots\left(
m-n+1\right)  }{n\left(  n-1\right)  !}=\dfrac{m\left(  m-1\right)
\cdots\left(  m-n+1\right)  }{n!}%
\end{align*}
(since $m\left(  m-1\right)  \left(  m-2\right)  \cdots\left(  m-n+1\right)
=m\left(  m-1\right)  \cdots\left(  m-n+1\right)  $ and $n\left(  n-1\right)
!=n!$). Compared with (\ref{eq.binom.mn}), this yields $\dbinom{m}{n}%
=\dfrac{m}{n}\dbinom{m-1}{n-1}$. This proves Proposition \ref{prop.binom.X-1}.
\end{proof}

The relation (\ref{eq.binom.X-1}) is called the \textit{absorption identity}
in \cite[\S 5.1]{GKP}.

\begin{exercise}
\label{exe.binom.scaryfrac}Let $k$, $a$ and $b$ be three positive integers
such that $k\leq a\leq b$. Prove that%
\[
\dfrac{k-1}{k}\sum_{n=a}^{b}\dfrac{1}{\dbinom{n}{k}}=\dfrac{1}{\dbinom
{a-1}{k-1}}-\dfrac{1}{\dbinom{b}{k-1}}.
\]
(In particular, all fractions appearing in this equality are well-defined.)
\end{exercise}

\subsubsection{Trinomial revision}

\begin{proposition}
\label{prop.binom.trinom-rev}Let $m\in\mathbb{Q}$, $a\in\mathbb{N}$ and
$i\in\mathbb{N}$ be such that $i\geq a$. Then,%
\begin{equation}
\dbinom{m}{i}\dbinom{i}{a}=\dbinom{m}{a}\dbinom{m-a}{i-a}.
\label{eq.binom.trinom-rev.m}%
\end{equation}

\end{proposition}

\begin{proof}
[Proof of Proposition \ref{prop.binom.trinom-rev}.]Let $g=m-a$. Then, $m-a=g$.
Therefore,%
\begin{align*}
&  \dbinom{m}{a}\dbinom{m-a}{i-a}\\
&  =\underbrace{\dbinom{m}{a}}_{\substack{=\dfrac{m\left(  m-1\right)
\cdots\left(  m-a+1\right)  }{a!}\\\text{(by the definition of }\dbinom{m}%
{a}\text{)}}}\ \ \underbrace{\dbinom{g}{i-a}}_{\substack{=\dfrac{g\left(
g-1\right)  \cdots\left(  g-\left(  i-a\right)  +1\right)  }{\left(
i-a\right)  !}\\\text{(by the definition of }\dbinom{g}{i-a}\text{)}}}\\
&  =\dfrac{m\left(  m-1\right)  \cdots\left(  m-a+1\right)  }{a!}\cdot
\dfrac{g\left(  g-1\right)  \cdots\left(  g-\left(  i-a\right)  +1\right)
}{\left(  i-a\right)  !}\\
&  =\dfrac{m\left(  m-1\right)  \cdots\left(  g+1\right)  }{a!}\cdot
\dfrac{g\left(  g-1\right)  \cdots\left(  m-i+1\right)  }{\left(  i-a\right)
!}\\
&  \ \ \ \ \ \ \ \ \ \ \ \ \ \ \ \ \ \ \ \ \left(  \text{since }m-a=g\text{
and }\underbrace{g}_{=m-a}-\left(  i-a\right)  =\left(  m-a\right)  -\left(
i-a\right)  =m-i\right) \\
&  =\dfrac{1}{a!\cdot\left(  i-a\right)  !}\cdot\underbrace{\left(  m\left(
m-1\right)  \cdots\left(  g+1\right)  \right)  \cdot\left(  g\left(
g-1\right)  \cdots\left(  m-i+1\right)  \right)  }_{\substack{=m\left(
m-1\right)  \cdots\left(  m-i+1\right)  }}\\
&  =\dfrac{1}{a!\cdot\left(  i-a\right)  !}\cdot m\left(  m-1\right)
\cdots\left(  m-i+1\right)  .
\end{align*}
Compared with%
\begin{align*}
&  \underbrace{\dbinom{m}{i}}_{\substack{=\dfrac{m\left(  m-1\right)
\cdots\left(  m-i+1\right)  }{i!}\\\text{(by the definition of }\dbinom{m}%
{i}\text{)}}}\ \ \underbrace{\dbinom{i}{a}}_{\substack{=\dfrac{i!}{a!\left(
i-a\right)  !}\\\text{(by (\ref{eq.binom.formula}), applied to }i\text{ and
}a\\\text{instead of }m\text{ and }n\text{)}}}\\
&  =\dfrac{m\left(  m-1\right)  \cdots\left(  m-i+1\right)  }{i!}\cdot
\dfrac{i!}{a!\left(  i-a\right)  !}=\dfrac{1}{a!\cdot\left(  i-a\right)
!}\cdot m\left(  m-1\right)  \cdots\left(  m-i+1\right)  ,
\end{align*}
this yields $\dbinom{m}{i}\dbinom{i}{a}=\dbinom{m}{a}\dbinom{m-a}{i-a}$. This
proves Proposition \ref{prop.binom.trinom-rev}.

[Notice that we used (\ref{eq.binom.formula}) to simplify $\dbinom{i}{a}$ in
this proof. Do not be tempted to use (\ref{eq.binom.formula}) to simplify
$\dbinom{m}{i}$, $\dbinom{m}{a}$ and $\dbinom{m-a}{i-a}$: The $m$ in these
expressions may fail to be an integer!]
\end{proof}

Proposition \ref{prop.binom.trinom-rev} is a simple and yet highly useful
formula, which Graham, Knuth and Patashnik call \textit{trinomial revision} in
\cite[Table 174]{GKP}.

\subsection{Binomial coefficients and polynomials}

We have so far defined the binomial coefficient $\dbinom{m}{n}$ in the case
when $n\in\mathbb{N}$ while $m$ is some number (rational, real or complex).
However, we can take this definition even further: For example, we can define
$\dbinom{m}{n}$ when $m$ is a polynomial with rational or real coefficients.
Let us do this now:

\begin{definition}
Let $n\in\mathbb{N}$. Let $m$ be a polynomial whose coefficients are rational
numbers (or real numbers, or complex numbers).

We define a polynomial $\dbinom{m}{n}$ by the equality (\ref{eq.binom.mn}).
This is a polynomial whose coefficients will be rational numbers or real
numbers or complex numbers, depending on the coefficients of $m$.
\end{definition}

Thus, in particular, for the polynomial $X\in\mathbb{Q}\left[  X\right]  $, we
have%
\[
\dbinom{X}{n}=\dfrac{X\left(  X-1\right)  \cdots\left(  X-n+1\right)  }%
{n!}\ \ \ \ \ \ \ \ \ \ \text{for every }n\in\mathbb{N}.
\]
In particular,%
\begin{align*}
\dbinom{X}{0}  &  =\dfrac{X\left(  X-1\right)  \cdots\left(  X-0+1\right)
}{0!}=\dfrac{\left(  \text{empty product}\right)  }{1}=1;\\
\dbinom{X}{1}  &  =\dfrac{X\left(  X-1\right)  \cdots\left(  X-1+1\right)
}{1!}=\dfrac{X}{1}=X;\\
\dbinom{X}{2}  &  =\dfrac{X\left(  X-1\right)  }{2!}=\dfrac{X\left(
X-1\right)  }{2}=\dfrac{1}{2}X^{2}-\dfrac{1}{2}X;\\
\dbinom{X}{3}  &  =\dfrac{X\left(  X-1\right)  \left(  X-2\right)  }%
{3!}=\dfrac{X\left(  X-1\right)  \left(  X-2\right)  }{6}=\dfrac{1}{6}%
X^{3}-\dfrac{1}{2}X^{2}+\dfrac{1}{3}X.
\end{align*}

The polynomial $\dbinom{X}{n}$ lets us compute the binomial coefficients
$\dbinom{m}{n}$ for all $m\in\mathbb{N}$, because of the following:

\begin{proposition}
\label{prop.binom.mn}Let $m\in\mathbb{Q}$ and $n\in\mathbb{N}$. Then, the
rational number $\dbinom{m}{n}$ is the result of evaluating the polynomial
$\dbinom{X}{n}$ at $X=m$.
\end{proposition}

\begin{proof}
[Proof of Proposition \ref{prop.binom.mn}.]We have $\dbinom{X}{n}%
=\dfrac{X\left(  X-1\right)  \cdots\left(  X-n+1\right)  }{n!}$. Hence, the
result of evaluating the polynomial $\dbinom{X}{n}$ at $X=m$ is%
\[
\dfrac{m\left(  m-1\right)  \cdots\left(  m-n+1\right)  }{n!}=\dbinom{m}%
{n}\ \ \ \ \ \ \ \ \ \ \left(  \text{by (\ref{eq.binom.mn})}\right)  .
\]
This proves Proposition \ref{prop.binom.mn}.
\end{proof}

We note the following properties of the polynomials $\dbinom{X}{n}$:

\begin{proposition}
\label{prop.binom.Xes}\textbf{(a)} We have
\[
\dbinom{X}{0}=1.
\]

\textbf{(b)} We have%
\[
\dbinom{X}{1}=X.
\]

\textbf{(c)} For every $n\in\left\{  1,2,3,\ldots\right\}  $, we have%
\[
\dbinom{X}{n}=\dbinom{X-1}{n}+\dbinom{X-1}{n-1}.
\]

\textbf{(d)} For every $n\in\mathbb{N}$, we have%
\[
\dbinom{X}{n}=\left(  -1\right)  ^{n}\dbinom{n-X-1}{n}.
\]

\textbf{(e)} For every $n\in\left\{  1,2,3,\ldots\right\}  $, we have%
\[
\dbinom{X}{n}=\dfrac{X}{n}\dbinom{X-1}{n-1}.
\]

\textbf{(f)} Let $a\in\mathbb{N}$ and $i\in\mathbb{N}$ be such that $i\geq a$.
Then,%
\[
\dbinom{X}{i}\dbinom{i}{a}=\dbinom{X}{a}\dbinom{X-a}{i-a}.
\]

\end{proposition}

\begin{proof}
[Proof of Proposition \ref{prop.binom.Xes}.]\textbf{(a)} To obtain a proof of
Proposition \ref{prop.binom.Xes} \textbf{(a)}, replace every appearance of
\textquotedblleft$m$\textquotedblright\ by \textquotedblleft$X$%
\textquotedblright\ in the proof of Proposition \ref{prop.binom.00}
\textbf{(a)}.

\textbf{(b)} To obtain a proof of Proposition \ref{prop.binom.Xes}
\textbf{(b)}, replace every appearance of \textquotedblleft$m$%
\textquotedblright\ by \textquotedblleft$X$\textquotedblright\ in the proof of
Proposition \ref{prop.binom.00} \textbf{(b)}.

\textbf{(c)} To obtain a proof of Proposition \ref{prop.binom.Xes}
\textbf{(c)}, replace every appearance of \textquotedblleft$m$%
\textquotedblright\ by \textquotedblleft$X$\textquotedblright\ in the proof of
Proposition \ref{prop.binom.rec}.

\textbf{(d)} To obtain a proof of Proposition \ref{prop.binom.Xes}
\textbf{(d)}, replace every appearance of \textquotedblleft$m$%
\textquotedblright\ by \textquotedblleft$X$\textquotedblright\ in the proof of
Proposition \ref{prop.binom.upper-neg}.

\textbf{(e)} To obtain a proof of Proposition \ref{prop.binom.Xes}
\textbf{(e)}, replace every appearance of \textquotedblleft$m$%
\textquotedblright\ by \textquotedblleft$X$\textquotedblright\ in the proof of
Proposition \ref{prop.binom.X-1}.

\textbf{(f)} To obtain a proof of Proposition \ref{prop.binom.Xes}
\textbf{(f)}, replace every appearance of \textquotedblleft$m$%
\textquotedblright\ by \textquotedblleft$X$\textquotedblright\ in the proof of
Proposition \ref{prop.binom.trinom-rev}.
\end{proof}

Recall that any polynomial $P\in\mathbb{Q}\left[  X\right]  $ (that is, any
polynomial in the indeterminate $X$ with rational coefficients) can be
quasi-uniquely written in the form $P\left(  X\right)  =\sum_{i=0}^{d}%
c_{i}X^{i}$ with rational $c_{0},c_{1},\ldots,c_{d}$. The word
\textquotedblleft quasi-uniquely\textquotedblright\ here means that the
coefficients $c_{0},c_{1},\ldots,c_{d}$ are uniquely determined when
$d\in\mathbb{N}$ is specified; they are not literally unique because we can
always increase $d$ by adding new $0$ coefficients (for example, the
polynomial $\left(  1+X\right)  ^{2}$ can be written both as $1+2X+X^{2}$ and
as $1+2X+X^{2}+0X^{3}+0X^{4}$).

It is not hard to check that an analogue of this statement holds with the
$X^{i}$ replaced by the $\dbinom{X}{i}$:

\begin{proposition}
\label{prop.hartshorne}\textbf{(a)} Any polynomial $P\in\mathbb{Q}\left[
X\right]  $ can be quasi-uniquely written in the form $P\left(  X\right)
=\sum_{i=0}^{d}c_{i}\dbinom{X}{i}$ with rational $c_{0},c_{1},\ldots,c_{d}$.
(Again, \textquotedblleft quasi-uniquely\textquotedblright\ means that we can
always increase $d$ by adding new $0$ coefficients, but apart from this the
$c_{0},c_{1},\ldots,c_{d}$ are uniquely determined.)

\textbf{(b)} The polynomial $P$ is \textit{integer-valued} (i.e., its values
at integers are integers) if and only if these rationals $c_{0},c_{1}%
,\ldots,c_{d}$ are integers.
\end{proposition}

We will not use this fact below, but it gives context to Theorem
\ref{thm.vandermonde.XY} and Exercise \ref{exe.ps1.1.2} further below. The
\textquotedblleft if\textquotedblright\ part of Proposition
\ref{prop.hartshorne} \textbf{(b)} follows from (\ref{eq.binom.int}). For a
full proof of Proposition \ref{prop.hartshorne} \textbf{(b)}, see
\cite[Theorem 10.3]{AndDosS}. See also \cite{daSilv12} for a proof of the
\textquotedblleft only if\textquotedblright\ part.

We shall now prove some facts and give some exercises about binomial
coefficients; but let us first prove a fundamental property of polynomials:

\begin{lemma}
\label{lem.polyeq}\textbf{(a)} Let $P$ be a polynomial in the indeterminate
$X$ with rational coefficients. Assume that $P\left(  x\right)  =0$ for all
$x\in\mathbb{N}$. Then, $P=0$ as polynomials\footnotemark.

\textbf{(b)} Let $P$ and $Q$ be two polynomials in the indeterminate $X$ with
rational coefficients. Assume that $P\left(  x\right)  =Q\left(  x\right)  $
for all $x\in\mathbb{N}$. Then, $P=Q$ as polynomials.

\textbf{(c)} Let $P$ be a polynomial in the indeterminates $X$ and $Y$ with
rational coefficients. Assume that $P\left(  x,y\right)  =0$ for all
$x\in\mathbb{N}$ and $y\in\mathbb{N}$. Then, $P=0$ as polynomials.

\textbf{(d)} Let $P$ and $Q$ be two polynomials in the indeterminates $X$ and
$Y$ with rational coefficients. Assume that $P\left(  x,y\right)  =Q\left(
x,y\right)  $ for all $x\in\mathbb{N}$ and $y\in\mathbb{N}$. Then, $P=Q$ as polynomials.
\end{lemma}

\footnotetext{Recall that two polynomials are said to be equal if and only if
their respective coefficients are equal.}

Lemma \ref{lem.polyeq} is a well-known property of polynomials with rational
coefficients; let us still prove it for the sake of completeness.

\begin{proof}
[Proof of Lemma \ref{lem.polyeq}.]\textbf{(a)} The polynomial $P$ satisfies
$P\left(  x\right)  =0$ for every $x\in\mathbb{N}$. Hence, every
$x\in\mathbb{N}$ is a root of $P$. Thus, the polynomial $P$ has infinitely
many roots. But a nonzero polynomial in one variable (with rational
coefficients) can only have finitely many roots\footnote{In fact, a stronger
statement holds: A nonzero polynomial in one variable (with rational
coefficients) having degree $n\geq0$ has at most $n$ roots. See, for example,
\cite[Corollary 1.8.24]{Goodman} or \cite[Theorem 1.58]{Joyce17} or
\cite[Corollary 4.5.10]{Walker87} or \cite[Corollary 33.7]{Elman18} or
\cite[Theorem 2.4.15]{Swanso18} or \cite[Corollary 1.14]{Knapp1} for a proof.
Note that Swanson, in \cite{Swanso18}, works with polynomial functions instead
of polynomials; but as far as roots are concerned, the difference does not
matter (since the roots of a polynomial are precisely the roots of the
corresponding polynomial function).}. If $P$ was nonzero, this would force a
contradiction with the sentence before. So $P$ must be zero. In other words,
$P=0$. Lemma \ref{lem.polyeq} \textbf{(a)} is proven.

\textbf{(b)} Every $x\in\mathbb{N}$ satisfies $\left(  P-Q\right)  \left(
x\right)  =P\left(  x\right)  -Q\left(  x\right)  =0$ (since $P\left(
x\right)  =Q\left(  x\right)  $). Hence, Lemma \ref{lem.polyeq} \textbf{(a)}
(applied to $P-Q$ instead of $P$) yields $P-Q=0$. Thus, $P=Q$. Lemma
\ref{lem.polyeq} \textbf{(b)} is thus proven.

\textbf{(c)} Every $x\in\mathbb{N}$ and $y\in\mathbb{N}$ satisfy%
\begin{equation}
P\left(  x,y\right)  =0. \label{pf.lem.polyeq.b.2}%
\end{equation}

We can write the polynomial $P$ in the form $P=\sum_{k=0}^{d}P_{k}\left(
X\right)  Y^{k}$, where $d$ is an integer and where each $P_{k}\left(
X\right)  $ (for $0\leq k\leq d$) is a polynomial in the single variable $X$.
Consider this $d$ and these $P_{k}\left(  X\right)  $.

Fix $\alpha\in\mathbb{N}$. Every $x\in\mathbb{N}$ satisfies%
\begin{align*}
P\left(  \alpha,x\right)   &  =\sum_{k=0}^{d}P_{k}\left(  \alpha\right)
x^{k}\\
&  \ \ \ \ \ \ \ \ \ \ \left(  \text{here, we substituted }\alpha\text{ and
}x\text{ for }X\text{ and }Y\text{ in }P=\sum_{k=0}^{d}P_{k}\left(  X\right)
Y^{k}\right)  ,
\end{align*}
so that $\sum_{k=0}^{d}P_{k}\left(  \alpha\right)  x^{k}=P\left(
\alpha,x\right)  =0$ (by (\ref{pf.lem.polyeq.b.2}), applied to $\alpha$ and
$x$ instead of $x$ and $y$).

Therefore, Lemma \ref{lem.polyeq} \textbf{(a)} (applied to $\sum_{k=0}%
^{d}P_{k}\left(  \alpha\right)  X^{k}$ instead of $P$) yields that
\[
\sum_{k=0}^{d}P_{k}\left(  \alpha\right)  X^{k}=0
\]
as polynomials (in the indeterminate $X$). In other words, all coefficients of
the polynomial $\sum_{k=0}^{d}P_{k}\left(  \alpha\right)  X^{k}$ are $0$. In
other words, $P_{k}\left(  \alpha\right)  =0$ for all $k\in\left\{
0,1,\ldots,d\right\}  $.

Now, let us forget that we fixed $\alpha$. We thus have shown that
$P_{k}\left(  \alpha\right)  =0$ for all $k\in\left\{  0,1,\ldots,d\right\}  $
and $\alpha\in\mathbb{N}$.

Let us now fix $k\in\left\{  0,1,\ldots,d\right\}  $. Then, $P_{k}\left(
\alpha\right)  =0$ for all $\alpha\in\mathbb{N}$. In other words,
$P_{k}\left(  x\right)  =0$ for all $x\in\mathbb{N}$. Hence, Lemma
\ref{lem.polyeq} \textbf{(a)} (applied to $P=P_{k}$) yields that $P_{k}=0$ as polynomials.

Let us forget that we fixed $k$. We thus have proven that $P_{k}=0$ as
polynomials for each $k\in\left\{  0,1,\ldots,d\right\}  $. Hence,
$P=\sum_{k=0}^{d}\underbrace{P_{k}\left(  X\right)  }_{=0}Y^{k}=0$. This
proves Lemma \ref{lem.polyeq} \textbf{(c)}.

\textbf{(d)} Every $x\in\mathbb{N}$ and $y\in\mathbb{N}$ satisfy%
\[
\left(  P-Q\right)  \left(  x,y\right)  =P\left(  x,y\right)  -Q\left(
x,y\right)  =0\ \ \ \ \ \ \ \ \ \ \left(  \text{since }P\left(  x,y\right)
=Q\left(  x,y\right)  \right)  .
\]
Hence, Lemma \ref{lem.polyeq} \textbf{(c)} (applied to $P-Q$ instead of $P$)
yields $P-Q=0$. Thus, $P=Q$. Lemma \ref{lem.polyeq} \textbf{(d)} is proven.
\end{proof}

Of course, Lemma \ref{lem.polyeq} can be generalized to polynomials in more
than two variables (the proof of Lemma \ref{lem.polyeq} \textbf{(c)}
essentially suggests how to prove this generalization by induction over the
number of variables).\footnote{If you know what a commutative ring is, you
might wonder whether Lemma \ref{lem.polyeq} can also be generalized to
polynomials with coefficients from other commutative rings (e.g., from
$\mathbb{R}$ or $\mathbb{C}$) instead of rational coefficients. In other
words, what happens if we replace \textquotedblleft rational
coefficients\textquotedblright\ by \textquotedblleft coefficients in
$R$\textquotedblright\ throughout Lemma \ref{lem.polyeq}, where $R$ is some
commutative ring? (Of course, we will then have to also replace $P\left(
x\right)  $ by $P\left(  x\cdot1_{R}\right)  $ and so on.)
\par
The answer is that Lemma \ref{lem.polyeq} becomes generally false if we don't
require anything more specific on $R$. However, there are certain conditions
on $R$ that make Lemma \ref{lem.polyeq} remain valid. For instance, Lemma
\ref{lem.polyeq} remains valid for $R=\mathbb{Z}$, for $R=\mathbb{R}$ and for
$R=\mathbb{C}$, as well as for $R$ being any polynomial ring over $\mathbb{Z}%
$, $\mathbb{Q}$, $\mathbb{R}$ or $\mathbb{C}$. More generally, Lemma
\ref{lem.polyeq} is valid if $R$ is any field of characteristic $0$ (i.e., any
field such that the elements $n\cdot1_{R}$ for $n$ ranging over $\mathbb{N}$
are pairwise distinct), or any subring of such a field.}

\subsection{The Chu-Vandermonde identity}

\subsubsection{The statements}

The following fact is known as the \textit{Chu-Vandermonde identity}%
\footnote{See
\href{https://en.wikipedia.org/wiki/Vandermonde\%27s_identity\%23Chu-Vandermonde_identity}{the
Wikipedia page} for part of its history. Usually, the equality $\dbinom
{x+y}{n}=\sum_{k=0}^{n}\dbinom{x}{k}\dbinom{y}{n-k}$ for two
\textbf{nonnegative integers} $x$ and $y$ (this is a particular case of
Theorem \ref{thm.vandermonde.rat}) is called the \textit{Vandermonde identity}
(or the \textit{Vandermonde convolution identity}), whereas the name
\textquotedblleft\textit{Chu-Vandermonde identity}\textquotedblright\ is used
for the identity $\dbinom{X+Y}{n}=\sum_{k=0}^{n}\dbinom{X}{k}\dbinom{Y}{n-k}$
in which $X$ and $Y$ are \textbf{indeterminates} (this is Theorem
\ref{thm.vandermonde.XY}). However, this seems to be mostly a matter of
convention (which isn't even universally followed); and anyway the two
identities are easily derived from one another as we will see in the second
proof of Theorem \ref{thm.vandermonde.XY}.}:

\begin{theorem}
\label{thm.vandermonde.rat}Let $n\in\mathbb{N}$, $x\in\mathbb{Q}$ and
$y\in\mathbb{Q}$. Then,%
\[
\dbinom{x+y}{n}=\sum_{k=0}^{n}\dbinom{x}{k}\dbinom{y}{n-k}.
\]

\end{theorem}

Let us also give an analogous statement for polynomials:

\begin{theorem}
\label{thm.vandermonde.XY}Let $n\in\mathbb{N}$. Then,%
\[
\dbinom{X+Y}{n}=\sum_{k=0}^{n}\dbinom{X}{k}\dbinom{Y}{n-k}%
\]
(an equality between polynomials in two variables $X$ and $Y$).
\end{theorem}

We will give two proofs of this theorem: one algebraic, and one
combinatorial.\footnote{Note that Theorem~\ref{thm.vandermonde.rat} appears in
\cite[(5.27)]{GKP}, where it is called \textit{Vandermonde's convolution}. The
second proof of Theorem~\ref{thm.vandermonde.rat} we shall show below is just
a more detailed writeup of the proof given there.}

\subsubsection{An algebraic proof}

\begin{proof}
[First proof of Theorem \ref{thm.vandermonde.rat}.]Forget that we fixed $n$,
$x$ and $y$. We thus must prove that for each $n\in\mathbb{N}$, we have%
\begin{equation}
\left(  \dbinom{x+y}{n}=\sum_{k=0}^{n}\dbinom{x}{k}\dbinom{y}{n-k}\text{ for
all }x\in\mathbb{Q}\text{ and }y\in\mathbb{Q}\right)  .
\label{pf.thm.vandermonde.rat.alg.goal}%
\end{equation}

We shall prove (\ref{pf.thm.vandermonde.rat.alg.goal}) by induction over $n$:

\textit{Induction base:} Let $x\in\mathbb{Q}$ and $y\in\mathbb{Q}$.
Proposition \ref{prop.binom.00} \textbf{(a)} (applied to $m=x$) yields
$\dbinom{x}{0}=1$. Proposition \ref{prop.binom.00} \textbf{(a)} (applied to
$m=y$) yields $\dbinom{y}{0}=1$. Now,%
\begin{equation}
\sum_{k=0}^{0}\dbinom{x}{k}\dbinom{y}{0-k}=\underbrace{\dbinom{x}{0}}%
_{=1}\underbrace{\dbinom{y}{0-0}}_{=\dbinom{y}{0}=1}=1.
\label{pf.thm.vandermonde.pf.2.1}%
\end{equation}

But Proposition \ref{prop.binom.00} \textbf{(a)} (applied to $m=x+y$) yields
$\dbinom{x+y}{0}=1$. Compared with (\ref{pf.thm.vandermonde.pf.2.1}), this
yields $\dbinom{x+y}{0}=\sum_{k=0}^{0}\dbinom{x}{k}\dbinom{y}{0-k}$.

Now, forget that we fixed $x$ and $y$. We thus have shown that%
\[
\dbinom{x+y}{0}=\sum_{k=0}^{0}\dbinom{x}{k}\dbinom{y}{0-k}\text{ for all }%
x\in\mathbb{Q}\text{ and }y\in\mathbb{Q}.
\]
In other words, (\ref{pf.thm.vandermonde.rat.alg.goal}) holds for $n=0$. This
completes the induction base.

\textit{Induction step:} Let $N$ be a positive integer. Assume that
(\ref{pf.thm.vandermonde.rat.alg.goal}) holds for $n=N-1$. We need to prove
that (\ref{pf.thm.vandermonde.rat.alg.goal}) holds for $n=N$. In other words,
we need to prove that%
\begin{equation}
\dbinom{x+y}{N}=\sum_{k=0}^{N}\dbinom{x}{k}\dbinom{y}{N-k}\text{ for all }%
x\in\mathbb{Q}\text{ and }y\in\mathbb{Q}. \label{pf.thm.vandermonde.pf.2.goal}%
\end{equation}

We have assumed that (\ref{pf.thm.vandermonde.rat.alg.goal}) holds for
$n=N-1$. In other words, we have%
\begin{equation}
\dbinom{x+y}{N-1}=\sum_{k=0}^{N-1}\dbinom{x}{k}\dbinom{y}{\left(  N-1\right)
-k}\text{ for all }x\in\mathbb{Q}\text{ and }y\in\mathbb{Q}.
\label{pf.thm.vandermonde.pf.2.4}%
\end{equation}

Now, let us prove (\ref{pf.thm.vandermonde.pf.2.goal}):

[\textit{Proof of} (\ref{pf.thm.vandermonde.pf.2.goal}): Fix $x\in\mathbb{Q}$
and $y\in\mathbb{Q}$. Then, (\ref{pf.thm.vandermonde.pf.2.4}) (applied to
$x-1$ instead of $x$) yields%
\begin{align*}
\dbinom{x-1+y}{N-1}  &  =\sum_{k=0}^{N-1}\dbinom{x-1}{k}\dbinom{y}{\left(
N-1\right)  -k}\\
&  =\sum_{k=1}^{N}\dbinom{x-1}{k-1}\underbrace{\dbinom{y}{\left(  N-1\right)
-\left(  k-1\right)  }}_{=\dbinom{y}{N-k}}\\
&  \ \ \ \ \ \ \ \ \ \ \left(  \text{here, we have substituted }k-1\text{ for
}k\text{ in the sum}\right) \\
&  =\sum_{k=1}^{N}\dbinom{x-1}{k-1}\dbinom{y}{N-k}.
\end{align*}
Since $x-1+y=x+y-1$, this rewrites as%
\begin{equation}
\dbinom{x+y-1}{N-1}=\sum_{k=1}^{N}\dbinom{x-1}{k-1}\dbinom{y}{N-k}.
\label{pf.thm.vandermonde.pf.2.5}%
\end{equation}
On the other hand, (\ref{pf.thm.vandermonde.pf.2.4}) (applied to $y-1$ instead
of $y$) shows that%
\begin{equation}
\dbinom{x+y-1}{N-1}=\sum_{k=0}^{N-1}\dbinom{x}{k}\underbrace{\dbinom
{y-1}{\left(  N-1\right)  -k}}_{=\dbinom{y-1}{N-k-1}}=\sum_{k=0}^{N-1}%
\dbinom{x}{k}\dbinom{y-1}{N-k-1}. \label{pf.thm.vandermonde.pf.2.6}%
\end{equation}

Next, we notice a simple consequence of (\ref{eq.binom.X-1}): We have%
\begin{equation}
\dfrac{x}{N}\dbinom{x-1}{a-1}=\dfrac{a}{N}\dbinom{x}{a}%
\ \ \ \ \ \ \ \ \ \ \text{for every }a\in\left\{  1,2,3,\ldots\right\}
\label{pf.thm.vandermonde.pf.2.7}%
\end{equation}
\footnote{\textit{Proof of (\ref{pf.thm.vandermonde.pf.2.7}):} Let
$a\in\left\{  1,2,3,\ldots\right\}  $. Then, (\ref{eq.binom.X-1}) (applied to
$m=x$ and $n=a$) yields $\dbinom{x}{a}=\dfrac{x}{a}\dbinom{x-1}{a-1}$. Hence,%
\[
\dfrac{a}{N}\underbrace{\dbinom{x}{a}}_{=\dfrac{x}{a}\dbinom{x-1}{a-1}%
}=\underbrace{\dfrac{a}{N}\cdot\dfrac{x}{a}}_{=\dfrac{x}{N}}\dbinom{x-1}%
{a-1}=\dfrac{x}{N}\dbinom{x-1}{a-1}.
\]
This proves (\ref{pf.thm.vandermonde.pf.2.7}).}. The same argument (applied to
$y$ instead of $x$) shows that%
\begin{equation}
\dfrac{y}{N}\dbinom{y-1}{a-1}=\dfrac{a}{N}\dbinom{y}{a}%
\ \ \ \ \ \ \ \ \ \ \text{for every }a\in\left\{  1,2,3,\ldots\right\}  .
\label{pf.thm.vandermonde.pf.2.8}%
\end{equation}

We have%
\begin{align*}
\dfrac{x}{N}\underbrace{\dbinom{x+y-1}{N-1}}_{\substack{=\sum_{k=1}^{N}%
\dbinom{x-1}{k-1}\dbinom{y}{N-k}\\\text{(by (\ref{pf.thm.vandermonde.pf.2.5}%
))}}}  &  =\dfrac{x}{N}\sum_{k=1}^{N}\dbinom{x-1}{k-1}\dbinom{y}{N-k}\\
&  =\sum_{k=1}^{N}\underbrace{\dfrac{x}{N}\dbinom{x-1}{k-1}}%
_{\substack{=\dfrac{k}{N}\dbinom{x}{k}\\\text{(by
(\ref{pf.thm.vandermonde.pf.2.7}),}\\\text{applied to }a=k\text{)}}}\dbinom
{y}{N-k}=\sum_{k=1}^{N}\dfrac{k}{N}\dbinom{x}{k}\dbinom{y}{N-k}.
\end{align*}
Compared with%
\begin{align*}
\sum_{k=0}^{N}\dfrac{k}{N}\dbinom{x}{k}\dbinom{y}{N-k}  &  =\underbrace{\dfrac
{0}{N}\dbinom{x}{0}\dbinom{y}{N-0}}_{=0}+\sum_{k=1}^{N}\dfrac{k}{N}\dbinom
{x}{k}\dbinom{y}{N-k}\\
&  \ \ \ \ \ \ \ \ \ \ \left(  \text{here, we have split off the addend for
}k=0\right) \\
&  =\sum_{k=1}^{N}\dfrac{k}{N}\dbinom{x}{k}\dbinom{y}{N-k},
\end{align*}
this yields%
\begin{equation}
\dfrac{x}{N}\dbinom{x+y-1}{N-1}=\sum_{k=0}^{N}\dfrac{k}{N}\dbinom{x}{k}%
\dbinom{y}{N-k}. \label{pf.thm.vandermonde.pf.2.11}%
\end{equation}

We also have%
\begin{align*}
\dfrac{y}{N}\underbrace{\dbinom{x+y-1}{N-1}}_{\substack{=\sum_{k=0}%
^{N-1}\dbinom{x}{k}\dbinom{y-1}{N-k-1}\\\text{(by
(\ref{pf.thm.vandermonde.pf.2.6}))}}}  &  =\dfrac{y}{N}\sum_{k=0}^{N-1}%
\dbinom{x}{k}\dbinom{y-1}{N-k-1}\\
&  =\sum_{k=0}^{N-1}\dbinom{x}{k}\underbrace{\dfrac{y}{N}\dbinom{y-1}{N-k-1}%
}_{\substack{=\dfrac{N-k}{N}\dbinom{y}{N-k}\\\text{(by
(\ref{pf.thm.vandermonde.pf.2.8}), applied to }a=N-k\\\text{(since }%
N-k\in\left\{  1,2,3,\ldots\right\}  \text{ (because }k\in\left\{
0,1,\ldots,N-1\right\}  \text{)))}}}\\
&  =\sum_{k=0}^{N-1}\dbinom{x}{k}\dfrac{N-k}{N}\dbinom{y}{N-k}=\sum
_{k=0}^{N-1}\dfrac{N-k}{N}\dbinom{x}{k}\dbinom{y}{N-k}.
\end{align*}
Compared with%
\begin{align*}
\sum_{k=0}^{N}\dfrac{N-k}{N}\dbinom{x}{k}\dbinom{y}{N-k}  &  =\sum_{k=0}%
^{N-1}\dfrac{N-k}{N}\dbinom{x}{k}\dbinom{y}{N-k}+\underbrace{\dfrac{N-N}{N}%
}_{=0}\dbinom{x}{N}\dbinom{y}{N-N}\\
&  \ \ \ \ \ \ \ \ \ \ \left(  \text{here, we have split off the addend for
}k=N\right) \\
&  =\sum_{k=0}^{N-1}\dfrac{N-k}{N}\dbinom{x}{k}\dbinom{y}{N-k}%
+\underbrace{0\dbinom{x}{N}\dbinom{y}{N-N}}_{=0}\\
&  =\sum_{k=0}^{N-1}\dfrac{N-k}{N}\dbinom{x}{k}\dbinom{y}{N-k},
\end{align*}
this yields%
\begin{equation}
\dfrac{y}{N}\dbinom{x+y-1}{N-1}=\sum_{k=0}^{N}\dfrac{N-k}{N}\dbinom{x}%
{k}\dbinom{y}{N-k}. \label{pf.thm.vandermonde.pf.2.12}%
\end{equation}

Now, (\ref{eq.binom.X-1}) (applied to $m=x+y$ and $n=N$) yields%
\begin{align*}
\dbinom{x+y}{N}  &  =\underbrace{\dfrac{x+y}{N}}_{=\dfrac{x}{N}+\dfrac{y}{N}%
}\dbinom{x+y-1}{N-1}=\left(  \dfrac{x}{N}+\dfrac{y}{N}\right)  \dbinom
{x+y-1}{N-1}\\
&  =\underbrace{\dfrac{x}{N}\dbinom{x+y-1}{N-1}}_{\substack{=\sum_{k=0}%
^{N}\dfrac{k}{N}\dbinom{x}{k}\dbinom{y}{N-k}\\\text{(by
(\ref{pf.thm.vandermonde.pf.2.11}))}}}+\underbrace{\dfrac{y}{N}\dbinom
{x+y-1}{N-1}}_{\substack{=\sum_{k=0}^{N}\dfrac{N-k}{N}\dbinom{x}{k}\dbinom
{y}{N-k}\\\text{(by (\ref{pf.thm.vandermonde.pf.2.12}))}}}\\
&  =\sum_{k=0}^{N}\dfrac{k}{N}\dbinom{x}{k}\dbinom{y}{N-k}+\sum_{k=0}%
^{N}\dfrac{N-k}{N}\dbinom{x}{k}\dbinom{y}{N-k}\\
&  =\sum_{k=0}^{N}\left(  \underbrace{\dfrac{k}{N}+\dfrac{N-k}{N}}%
_{=1}\right)  \dbinom{x}{k}\dbinom{y}{N-k}=\sum_{k=0}^{N}\dbinom{x}{k}%
\dbinom{y}{N-k}.
\end{align*}
This proves (\ref{pf.thm.vandermonde.pf.2.goal}).]

We have thus proven (\ref{pf.thm.vandermonde.pf.2.goal}). In other words,
(\ref{pf.thm.vandermonde.rat.alg.goal}) holds for $n=N$. This completes the
induction step. Thus, the induction proof of
(\ref{pf.thm.vandermonde.rat.alg.goal}) is complete. Hence, Theorem
\ref{thm.vandermonde.rat} is proven.
\end{proof}

The above proof has the advantage of being completely algebraic; it thus does
not rely on what $x$ and $y$ actually are. It works equally well if $x$ and
$y$ are assumed to be real numbers or complex numbers or polynomials. Thus, it
can also be used to prove Theorem \ref{thm.vandermonde.XY}:

\begin{proof}
[First proof of Theorem \ref{thm.vandermonde.XY}.]To obtain a proof of Theorem
\ref{thm.vandermonde.XY}, replace every appearance of \textquotedblleft%
$x$\textquotedblright\ by \textquotedblleft$X$\textquotedblright\ and every
appearance of \textquotedblleft$y$\textquotedblright\ by \textquotedblleft%
$Y$\textquotedblright\ in the above First proof of Theorem
\ref{thm.vandermonde.rat}.
\end{proof}

\subsubsection{A combinatorial proof}

We shall next give a different, combinatorial proof of Theorems
\ref{thm.vandermonde.rat} and \ref{thm.vandermonde.XY}. This proof is somewhat
indirect, as it begins by showing the following particular case of Theorem
\ref{thm.vandermonde.rat}:

\begin{lemma}
\label{lem.vandermonde.nat}Let $n\in\mathbb{N}$, $x\in\mathbb{N}$ and
$y\in\mathbb{N}$. Then,%
\begin{equation}
\dbinom{x+y}{n}=\sum_{k=0}^{n}\dbinom{x}{k}\dbinom{y}{n-k}.
\label{pf.thm.vandermonde.eq.xy}%
\end{equation}

\end{lemma}

This lemma is less general than Theorem \ref{thm.vandermonde.rat}, since it
requires $x$ and $y$ to belong to $\mathbb{N}$.

\begin{proof}
[Proof of Lemma \ref{lem.vandermonde.nat}.]For every $N\in\mathbb{N}$, we let
$\left[  N\right]  $ denote the $N$-element set $\left\{  1,2,\ldots
,N\right\}  $.

Recall that $\dbinom{x+y}{n}$ is the number of $n$-element subsets of a given
$\left(  x+y\right)  $-element set\footnote{This follows from
(\ref{eq.binom.subsets}).}. Since $\left[  x+y\right]  $ is an $\left(
x+y\right)  $-element set, we thus conclude that $\dbinom{x+y}{n}$ is the
number of $n$-element subsets of $\left[  x+y\right]  $.

But let us count the $n$-element subsets of $\left[  x+y\right]  $ in a
different way (i.e., find a different expression for the number of $n$-element
subsets of $\left[  x+y\right]  $). Namely, we can choose an $n$-element
subset $S$ of $\left[  x+y\right]  $ by means of the following process:

\begin{enumerate}
\item We decide how many elements of this subset $S$ will be among the numbers
$1,2,\ldots,x$. Let $k$ be the number of these elements. Clearly, $k$ must be
an integer between $0$ and $n$ (inclusive)\footnote{Because the subset $S$
will have $n$ elements in total, and thus at most $n$ of them can be among the
numbers $1,2,\ldots,x$.}.

\item Then, we choose these $k$ elements of $S$ among the numbers
$1,2,\ldots,x$. This can be done in $\dbinom{x}{k}$ different ways (because we
are choosing $k$ out of $x$ numbers, with no repetitions, and with no regard
for their order; in other words, we are choosing a $k$-element subset of
$\left\{  1,2,\ldots,x\right\}  $).

\item Then, we choose the remaining $n-k$ elements of $S$ (because $S$ should
have $n$ elements in total) among the remaining numbers $x+1,x+2,\ldots,x+y$.
This can be done in $\dbinom{y}{n-k}$ ways (because we are choosing $n-k$ out
of $y$ numbers, with no repetitions, and with no regard for their order).
\end{enumerate}

This process makes it clear that the total number of ways to choose an
$n$-element subset $S$ of $\left[  x+y\right]  $ is $\sum_{k=0}^{n}\dbinom
{x}{k}\dbinom{y}{n-k}$. In other words, the number of $n$-element subsets of
$\left[  x+y\right]  $ is $\sum_{k=0}^{n}\dbinom{x}{k}\dbinom{y}{n-k}$. But
earlier, we have shown that the same number is $\dbinom{x+y}{n}$. Comparing
these two results, we conclude that $\dbinom{x+y}{n}=\sum_{k=0}^{n}\dbinom
{x}{k}\dbinom{y}{n-k}$. Thus, Lemma \ref{lem.vandermonde.nat} is proven.
\end{proof}

We now shall leverage Lemma \ref{lem.polyeq} to derive Theorem
\ref{thm.vandermonde.XY} from this lemma:

\begin{proof}
[Second proof of Theorem \ref{thm.vandermonde.XY}.]We define two polynomials
$P$ and $Q$ in the indeterminates $X$ and $Y$ with rational coefficients by
setting%
\begin{align*}
P  &  =\dbinom{X+Y}{n};\\
Q  &  =\sum_{k=0}^{n}\dbinom{X}{k}\dbinom{Y}{n-k}%
\end{align*}
\footnote{These are both polynomials since $\dbinom{X+Y}{n}$, $\dbinom{X}{k}$
and $\dbinom{Y}{n-k}$ are polynomials in $X$ and $Y$.}. The equality
(\ref{pf.thm.vandermonde.eq.xy}) (which we have proven) states that $P\left(
x,y\right)  =Q\left(  x,y\right)  $ for all $x\in\mathbb{N}$ and
$y\in\mathbb{N}$. Thus, Lemma \ref{lem.polyeq} \textbf{(d)} yields that $P=Q$.
Recalling how $P$ and $Q$ are defined, we can rewrite this as $\dbinom{X+Y}%
{n}=\sum_{k=0}^{n}\dbinom{X}{k}\dbinom{Y}{n-k}$. This proves Theorem
\ref{thm.vandermonde.XY}.
\end{proof}

The argument that we used at the end of the above proof to derive Theorem
\ref{thm.vandermonde.XY} from (\ref{pf.thm.vandermonde.eq.xy}) is a very
common argument that appears in proofs of equalities for binomial
coefficients. The binomial coefficients $\dbinom{m}{n}$ are defined for
arbitrary rational, real or complex $m$\ \ \ \ \footnote{For example, terms
like $\dbinom{-1/2}{3}$, $\dbinom{2+\sqrt{3}}{5}$ and $\dbinom{-7}{0}$ make
perfect sense. (But we cannot substitute arbitrary complex numbers for $n$ in
$\dbinom{m}{n}$. So far we have only defined $\dbinom{m}{n}$ for
$n\in\mathbb{N}$. It is usual to define $\dbinom{m}{n}$ to mean $0$ for
negative integers $n$, and using analysis (specifically, the $\Gamma$
function) it is possible to give a reasonable meaning to $\dbinom{m}{n}$ for
$m$ and $n$ being reals, but this will no longer be a polynomial in $m$.)},
but their combinatorial interpretation (via counting subsets) only makes sense
when $m$ and $n$ are nonnegative integers. Thus, if we want to prove an
identity of the form $P=Q$ (where $P$ and $Q$ are two polynomials, say, in two
indeterminates $X$ and $Y$) using the combinatorial interpretation of binomial
coefficients, then a reasonable tactic is to first show that $P\left(
x,y\right)  =Q\left(  x,y\right)  $ for all $x\in\mathbb{N}$ and
$y\in\mathbb{N}$ (using combinatorics), and then to use something like Lemma
\ref{lem.polyeq} in order to conclude that $P$ and $Q$ are equal as
polynomials. We shall see this tactic used a few more times.\footnote{This
tactic is called \textquotedblleft the polynomial argument\textquotedblright%
\ in \cite[\S 5.1]{GKP}.}

\begin{proof}
[Second proof of Theorem \ref{thm.vandermonde.rat}.]Theorem
\ref{thm.vandermonde.XY} yields
\[
\dbinom{X+Y}{n}=\sum_{k=0}^{n}\dbinom{X}{k}\dbinom{Y}{n-k}%
\]
(an equality between polynomials in two variables $X$ and $Y$). Now, let us
evaluate both sides of this equality at $X=x$ and $Y=y$. As a result, we
obtain%
\[
\dbinom{x+y}{n}=\sum_{k=0}^{n}\dbinom{x}{k}\dbinom{y}{n-k}%
\]
(because of Proposition \ref{prop.binom.mn}).
\end{proof}

\subsubsection{Some applications}

Let us give some sample applications of Theorem \ref{thm.vandermonde.rat}:

\begin{proposition}
\label{prop.vandermonde.consequences}

\textbf{(a)} For every $x\in\mathbb{Z}$ and $y\in\mathbb{Z}$ and
$n\in\mathbb{N}$, we have
\[
\dbinom{x+y}{n}=\sum_{k=0}^{n}\dbinom{x}{k}\dbinom{y}{n-k}.
\]

\textbf{(b)} For every $x\in\mathbb{N}$ and $y\in\mathbb{Z}$, we have%
\[
\dbinom{x+y}{x}=\sum_{k=0}^{x}\dbinom{x}{k}\dbinom{y}{k}.
\]

\textbf{(c)} For every $n\in\mathbb{N}$, we have%
\[
\dbinom{2n}{n}=\sum_{k=0}^{n}\dbinom{n}{k}^{2}.
\]

\textbf{(d)} For every $x\in\mathbb{Z}$ and $y\in\mathbb{Z}$ and
$n\in\mathbb{N}$, we have%
\[
\dbinom{x-y}{n}=\sum_{k=0}^{n}\left(  -1\right)  ^{k}\dbinom{x}{n-k}%
\dbinom{k+y-1}{k}.
\]

\textbf{(e)} For every $x\in\mathbb{N}$ and $y\in\mathbb{Z}$ and
$n\in\mathbb{N}$ with $x\leq n$, we have%
\[
\dbinom{y-x-1}{n-x}=\sum_{k=0}^{n}\left(  -1\right)  ^{k-x}\dbinom{k}%
{x}\dbinom{y}{n-k}.
\]

\textbf{(f)} For every $x\in\mathbb{N}$ and $y\in\mathbb{N}$ and
$n\in\mathbb{N}$, we have%
\[
\dbinom{n+1}{x+y+1}=\sum_{k=0}^{n}\dbinom{k}{x}\dbinom{n-k}{y}.
\]

\textbf{(g)} For every $x\in\mathbb{Z}$ and $y\in\mathbb{N}$ and
$n\in\mathbb{N}$ satisfying $x+y\geq0$ and $n\geq x$, we have%
\[
\dbinom{x+y}{n}=\sum_{k=0}^{x+y}\dbinom{x}{k}\dbinom{y}{n+k-x}.
\]

\end{proposition}

\begin{remark}
I have learnt Proposition \ref{prop.vandermonde.consequences} \textbf{(f)}
from \href{http://www.artofproblemsolving.com/community/c6h447278}{the AoPS
forum}. Proposition \ref{prop.vandermonde.consequences} \textbf{(g)} is a
generalization of Proposition \ref{prop.vandermonde.consequences} \textbf{(b)}.

Note that if we apply Proposition \ref{prop.vandermonde.consequences}
\textbf{(f)} to $y=0$, then we obtain the identity $\dbinom{n+1}{x+1}%
=\sum_{k=0}^{n}\dbinom{k}{x}$ for all $n\in\mathbb{N}$ and $x\in\mathbb{N}$.
This identity is also a particular case of Exercise \ref{exe.binom.hockey1}
\textbf{(b)}.
\end{remark}

\begin{proof}
[Proof of Proposition \ref{prop.vandermonde.consequences}.]\textbf{(a)}
Proposition \ref{prop.vandermonde.consequences} \textbf{(a)} is a particular
case of Theorem \ref{thm.vandermonde.rat}.

\textbf{(b)} Let $x\in\mathbb{N}$ and $y\in\mathbb{Z}$. Proposition
\ref{prop.vandermonde.consequences} \textbf{(a)} (applied to $y$, $x$ and $x$
instead of $x$, $y$ and $n$) yields%
\[
\dbinom{y+x}{x}=\sum_{k=0}^{x}\dbinom{y}{k}\dbinom{x}{x-k}.
\]
Compared with
\[
\sum_{k=0}^{x}\underbrace{\dbinom{x}{k}}_{\substack{=\dbinom{x}{x-k}%
\\\text{(by (\ref{eq.binom.symm}), applied to }m=x\text{ and }%
n=k\\\text{(since }x\geq k\text{ (because }k\leq x\text{)))}}}\dbinom{y}%
{k}=\sum_{k=0}^{x}\dbinom{x}{x-k}\dbinom{y}{k}=\sum_{k=0}^{x}\dbinom{y}%
{k}\dbinom{x}{x-k},
\]
this yields $\dbinom{y+x}{x}=\sum_{k=0}^{x}\dbinom{x}{k}\dbinom{y}{k}$. Since
$y+x=x+y$, this rewrites as $\dbinom{x+y}{x}=\sum_{k=0}^{x}\dbinom{x}%
{k}\dbinom{y}{k}$. This proves Proposition \ref{prop.vandermonde.consequences}
\textbf{(b)}.

\textbf{(c)} Let $n\in\mathbb{N}$. Applying Proposition
\ref{prop.vandermonde.consequences} \textbf{(b)} to $x=n$ and $y=n$, we obtain%
\[
\dbinom{n+n}{n}=\sum_{k=0}^{n}\underbrace{\dbinom{n}{k}\dbinom{n}{k}%
}_{=\dbinom{n}{k}^{2}}=\sum_{k=0}^{n}\dbinom{n}{k}^{2}.
\]
Since $n+n=2n$, this rewrites as $\dbinom{2n}{n}=\sum_{k=0}^{n}\dbinom{n}%
{k}^{2}$. This proves Proposition \ref{prop.vandermonde.consequences}
\textbf{(c)}.

\textbf{(d)} Let $x\in\mathbb{Z}$ and $y\in\mathbb{Z}$ and $n\in\mathbb{N}$.
Proposition \ref{prop.vandermonde.consequences} \textbf{(a)} (applied to $-y$
instead of $y$) yields%
\begin{align*}
\dbinom{x+\left(  -y\right)  }{n}  &  =\sum_{k=0}^{n}\dbinom{x}{k}\dbinom
{-y}{n-k}=\sum_{k=0}^{n}\dbinom{x}{n-k}\underbrace{\dbinom{-y}{n-\left(
n-k\right)  }}_{=\dbinom{-y}{k}}\\
&  \ \ \ \ \ \ \ \ \ \ \left(  \text{here, we substituted }n-k\text{ for
}k\text{ in the sum}\right) \\
&  =\sum_{k=0}^{n}\dbinom{x}{n-k}\underbrace{\dbinom{-y}{k}}%
_{\substack{=\left(  -1\right)  ^{k}\dbinom{k-\left(  -y\right)  -1}%
{k}\\\text{(by (\ref{eq.binom.upper-neg}), applied to }k\text{ and }-y\text{
instead of }n\text{ and }m\text{)}}}\\
&  =\sum_{k=0}^{n}\underbrace{\dbinom{x}{n-k}\left(  -1\right)  ^{k}%
}_{=\left(  -1\right)  ^{k}\dbinom{x}{n-k}}\underbrace{\dbinom{k-\left(
-y\right)  -1}{k}}_{=\dbinom{k+y-1}{k}}\\
&  =\sum_{k=0}^{n}\left(  -1\right)  ^{k}\dbinom{x}{n-k}\dbinom{k+y-1}{k}.
\end{align*}
Since $x+\left(  -y\right)  =x-y$, this rewrites as $\dbinom{x-y}{n}%
=\sum_{k=0}^{n}\left(  -1\right)  ^{k}\dbinom{x}{n-k}\dbinom{k+y-1}{k}$. This
proves Proposition \ref{prop.vandermonde.consequences} \textbf{(d)}.

\textbf{(e)} Let $x\in\mathbb{N}$ and $y\in\mathbb{Z}$ and $n\in\mathbb{N}$ be
such that $x\leq n$. From $x\in\mathbb{N}$, we obtain $0\leq x$ and thus
$0\leq x\leq n$. We notice that every integer $k\geq x$ satisfies%
\begin{equation}
\dbinom{k}{k-x}=\dbinom{k}{x} \label{pf.prop.vandermonde.consequences.e.1}%
\end{equation}
\footnote{\textit{Proof of (\ref{pf.prop.vandermonde.consequences.e.1}):} Let
$k$ be an integer such that $k\geq x$. Thus, $k-x\in\mathbb{N}$. Also, $k\geq
x\geq0$ (since $x\in\mathbb{N}$), and thus $k\in\mathbb{N}$. Now, recall that
$k\geq x$. Hence, (\ref{eq.binom.symm}) (applied to $k$ and $x$ instead of $m$
and $n$) yields $\dbinom{k}{x}=\dbinom{k}{k-x}$. This proves
(\ref{pf.prop.vandermonde.consequences.e.1}).}. Furthermore, $n-x\in
\mathbb{N}$ (since $x\leq n$). Hence, we can apply Proposition
\ref{prop.vandermonde.consequences} \textbf{(a)} to $y$, $-x-1$ and $n-x$
instead of $x$, $y$ and $n$. As a result, we obtain%
\begin{align*}
\dbinom{y+\left(  -x-1\right)  }{n-x}  &  =\sum_{k=0}^{n-x}\dbinom{y}%
{k}\underbrace{\dbinom{-x-1}{\left(  n-x\right)  -k}}_{\substack{=\left(
-1\right)  ^{\left(  n-x\right)  -k}\dbinom{\left(  \left(  n-x\right)
-k\right)  -\left(  -x-1\right)  -1}{\left(  n-x\right)  -k}\\\text{(by
(\ref{eq.binom.upper-neg}), applied to }-x-1\text{ and }\left(  n-x\right)
-k\\\text{instead of }m\text{ and }n\text{)}}}\\
&  =\sum_{k=0}^{n-x}\dbinom{y}{k}\left(  -1\right)  ^{\left(  n-x\right)
-k}\underbrace{\dbinom{\left(  \left(  n-x\right)  -k\right)  -\left(
-x-1\right)  -1}{\left(  n-x\right)  -k}}_{\substack{=\dbinom{n-k}{\left(
n-x\right)  -k}\\\text{(since }\left(  \left(  n-x\right)  -k\right)  -\left(
-x-1\right)  -1=n-k\text{)}}}\\
&  =\sum_{k=0}^{n-x}\dbinom{y}{k}\left(  -1\right)  ^{\left(  n-x\right)
-k}\dbinom{n-k}{\left(  n-x\right)  -k}\\
&  =\underbrace{\sum_{k=n-\left(  n-x\right)  }^{n}}_{\substack{=\sum
_{k=x}^{n}\\\text{(since }n-\left(  n-x\right)  =x\text{)}}}\dbinom{y}%
{n-k}\underbrace{\left(  -1\right)  ^{\left(  n-x\right)  -\left(  n-k\right)
}}_{\substack{=\left(  -1\right)  ^{k-x}\\\text{(since }\left(  n-x\right)
-\left(  n-k\right)  =k-x\text{)}}}\underbrace{\dbinom{n-\left(  n-k\right)
}{\left(  n-x\right)  -\left(  n-k\right)  }}_{\substack{=\dbinom{k}%
{k-x}\\\text{(since }n-\left(  n-k\right)  =k\\\text{and }\left(  n-x\right)
-\left(  n-k\right)  =k-x\text{)}}}\\
&  \ \ \ \ \ \ \ \ \ \ \left(  \text{here, we have substituted }n-k\text{ for
}k\text{ in the sum}\right) \\
&  =\sum_{k=x}^{n}\dbinom{y}{n-k}\left(  -1\right)  ^{k-x}\underbrace{\dbinom
{k}{k-x}}_{\substack{=\dbinom{k}{x}\\\text{(by
(\ref{pf.prop.vandermonde.consequences.e.1}))}}}=\sum_{k=x}^{n}\dbinom{y}%
{n-k}\left(  -1\right)  ^{k-x}\dbinom{k}{x}\\
&  =\sum_{k=x}^{n}\left(  -1\right)  ^{k-x}\dbinom{k}{x}\dbinom{y}{n-k}.
\end{align*}
Compared with%
\begin{align*}
&  \sum_{k=0}^{n}\left(  -1\right)  ^{k-x}\dbinom{k}{x}\dbinom{y}{n-k}\\
&  =\sum_{k=0}^{x-1}\left(  -1\right)  ^{k-x}\underbrace{\dbinom{k}{x}%
}_{\substack{=0\\\text{(by (\ref{eq.binom.0}), applied to }k\text{ and
}x\\\text{instead of }m\text{ and }n\text{ (since }k\leq x-1<x\text{))}%
}}\dbinom{y}{n-k}+\sum_{k=x}^{n}\left(  -1\right)  ^{k-x}\dbinom{k}{x}%
\dbinom{y}{n-k}\\
&  \ \ \ \ \ \ \ \ \ \ \left(  \text{since }0\leq x\leq n\right) \\
&  =\underbrace{\sum_{k=0}^{x-1}\left(  -1\right)  ^{k-x}0\dbinom{y}{n-k}%
}_{=0}+\sum_{k=x}^{n}\left(  -1\right)  ^{k-x}\dbinom{k}{x}\dbinom{y}%
{n-k}=\sum_{k=x}^{n}\left(  -1\right)  ^{k-x}\dbinom{k}{x}\dbinom{y}{n-k},
\end{align*}
this yields%
\[
\dbinom{y+\left(  -x-1\right)  }{n-x}=\sum_{k=0}^{n}\left(  -1\right)
^{k-x}\dbinom{k}{x}\dbinom{y}{n-k}.
\]
In other words,
\[
\dbinom{y-x-1}{n-x}=\sum_{k=0}^{n}\left(  -1\right)  ^{k-x}\dbinom{k}%
{x}\dbinom{y}{n-k}%
\]
(since $y+\left(  -x-1\right)  =y-x-1$). This proves Proposition
\ref{prop.vandermonde.consequences} \textbf{(e)}.

\textbf{(f)} Let $x\in\mathbb{N}$ and $y\in\mathbb{N}$ and $n\in\mathbb{N}$.
We must be in one of the following two cases:

\textit{Case 1:} We have $n<x+y$.

\textit{Case 2:} We have $n\geq x+y$.

Let us first consider Case 1. In this case, we have $n<x+y$. Thus,
$n+1<x+y+1$. Therefore, $\dbinom{n+1}{x+y+1}=0$ (by (\ref{eq.binom.0}),
applied to $n+1$ and $x+y+1$ instead of $m$ and $n$). But every $k\in\left\{
0,1,\ldots,n\right\}  $ satisfies $\dbinom{k}{x}\dbinom{n-k}{y}=0$%
\ \ \ \ \footnote{\textit{Proof.} Let $k\in\left\{  0,1,\ldots,n\right\}  $.
We need to show that $\dbinom{k}{x}\dbinom{n-k}{y}=0$.
\par
If we have $k<x$, then we have $\dbinom{k}{x}=0$ (by (\ref{eq.binom.0}),
applied to $k$ and $x$ instead of $m$ and $n$). Therefore, if we have $k<x$,
then $\underbrace{\dbinom{k}{x}}_{=0}\dbinom{n-k}{y}=0$. Hence, for the rest
of this proof of $\dbinom{k}{x}\dbinom{n-k}{y}=0$, we can WLOG assume that we
don't have $k<x$. Assume this.
\par
We have $k\leq n$ (since $k\in\left\{  0,1,\ldots,n\right\}  $) and thus
$n-k\in\mathbb{N}$.
\par
We have $k\geq x$ (since we don't have $k<x$), and thus $n-\underbrace{k}%
_{\geq x}\leq n-x<y$ (since $n<x+y$). Hence, $\dbinom{n-k}{y}=0$ (by
(\ref{eq.binom.0}), applied to $n-k$ and $y$ instead of $m$ and $n$).
Therefore, $\dbinom{k}{x}\underbrace{\dbinom{n-k}{y}}_{=0}=0$, qed.}. Hence,
$\sum_{k=0}^{n}\underbrace{\dbinom{k}{x}\dbinom{n-k}{y}}_{=0}=\sum_{k=0}%
^{n}0=0$. Compared with $\dbinom{n+1}{x+y+1}=0$, this yields $\dbinom
{n+1}{x+y+1}=\sum_{k=0}^{n}\dbinom{k}{x}\dbinom{n-k}{y}$. Thus, Proposition
\ref{prop.vandermonde.consequences} \textbf{(f)} is proven in Case 1.

Let us now consider Case 2. In this case, we have $n\geq x+y$. Hence, $n-y\geq
x$ (since $x\in\mathbb{N}$), so that $\left(  n-y\right)  -x\in\mathbb{N}$.
Also, $n-y\geq x\geq0$ and thus $n-y\in\mathbb{N}$. Moreover, $x\leq n-y$.
Therefore, we can apply Proposition \ref{prop.vandermonde.consequences}
\textbf{(e)} to $-y-1$ and $n-y$ instead of $y$ and $n$. As a result, we
obtain%
\begin{align}
&  \dbinom{\left(  -y-1\right)  -x-1}{\left(  n-y\right)  -x}\nonumber\\
&  =\sum_{k=0}^{n-y}\left(  -1\right)  ^{k-x}\dbinom{k}{x}\underbrace{\dbinom
{-y-1}{\left(  n-y\right)  -k}}_{\substack{=\left(  -1\right)  ^{\left(
n-y\right)  -k}\dbinom{\left(  \left(  n-y\right)  -k\right)  -\left(
-y-1\right)  -1}{\left(  n-y\right)  -k}\\\text{(by (\ref{eq.binom.upper-neg}%
), applied to }-y-1\text{ and }\left(  n-y\right)  -k\\\text{instead of
}m\text{ and }n\text{)}}}\nonumber\\
&  =\sum_{k=0}^{n-y}\left(  -1\right)  ^{k-x}\underbrace{\dbinom{k}{x}\left(
-1\right)  ^{\left(  n-y\right)  -k}}_{=\left(  -1\right)  ^{\left(
n-y\right)  -k}\dbinom{k}{x}}\underbrace{\dbinom{\left(  \left(  n-y\right)
-k\right)  -\left(  -y-1\right)  -1}{\left(  n-y\right)  -k}}%
_{\substack{=\dbinom{n-k}{\left(  n-y\right)  -k}\\\text{(since }\left(
\left(  n-y\right)  -k\right)  -\left(  -y-1\right)  -1=n-k\text{)}%
}}\nonumber\\
&  =\sum_{k=0}^{n-y}\underbrace{\left(  -1\right)  ^{k-x}\left(  -1\right)
^{\left(  n-y\right)  -k}}_{\substack{=\left(  -1\right)  ^{\left(
k-x\right)  +\left(  \left(  n-y\right)  -k\right)  }=\left(  -1\right)
^{n-x-y}\\\text{(since }\left(  k-x\right)  +\left(  \left(  n-y\right)
-k\right)  =n-x-y\text{)}}}\dbinom{k}{x}\dbinom{n-k}{\left(  n-y\right)
-k}\nonumber\\
&  =\sum_{k=0}^{n-y}\left(  -1\right)  ^{n-x-y}\dbinom{k}{x}\dbinom
{n-k}{\left(  n-y\right)  -k}\nonumber\\
&  =\left(  -1\right)  ^{n-x-y}\sum_{k=0}^{n-y}\dbinom{k}{x}\dbinom
{n-k}{\left(  n-y\right)  -k}. \label{pf.prop.vandermonde.consequences.f.0}%
\end{align}

But every $k\in\left\{  0,1,\ldots,n-y\right\}  $ satisfies%
\begin{equation}
\dbinom{n-k}{\left(  n-y\right)  -k}=\dbinom{n-k}{y}
\label{pf.prop.vandermonde.consequences.f.1}%
\end{equation}
\footnote{\textit{Proof of (\ref{pf.prop.vandermonde.consequences.f.1}):} Let
$k\in\left\{  0,1,\ldots,n-y\right\}  $. Then, $k\in\mathbb{N}$ and $n-y\geq
k$. From $n-y\geq k$, we obtain $n\geq y+k$, so that $n-k\geq y$. Thus,
$n-k\geq y\geq0$, so that $n-k\in\mathbb{N}$. Hence, (\ref{eq.binom.symm})
(applied to $n-k$ and $y$ instead of $m$ and $n$) yields $\dbinom{n-k}%
{y}=\dbinom{n-k}{\left(  n-k\right)  -y}=\dbinom{n-k}{\left(  n-y\right)  -k}$
(since $\left(  n-k\right)  -y=\left(  n-y\right)  -k$). This proves
(\ref{pf.prop.vandermonde.consequences.f.1}).}. Thus,
(\ref{pf.prop.vandermonde.consequences.f.0}) yields%
\begin{align*}
\dbinom{\left(  -y-1\right)  -x-1}{\left(  n-y\right)  -x}  &  =\left(
-1\right)  ^{n-x-y}\sum_{k=0}^{n-y}\dbinom{k}{x}\underbrace{\dbinom
{n-k}{\left(  n-y\right)  -k}}_{\substack{=\dbinom{n-k}{y}\\\text{(by
(\ref{pf.prop.vandermonde.consequences.f.1}))}}}\\
&  =\left(  -1\right)  ^{n-x-y}\sum_{k=0}^{n-y}\dbinom{k}{x}\dbinom{n-k}{y}.
\end{align*}
Compared with%
\begin{align*}
\dbinom{\left(  -y-1\right)  -x-1}{\left(  n-y\right)  -x}  &
=\underbrace{\left(  -1\right)  ^{\left(  n-y\right)  -x}}_{\substack{=\left(
-1\right)  ^{n-x-y}\\\text{(since }\left(  n-y\right)  -x=n-x-y\text{)}%
}}\underbrace{\dbinom{\left(  \left(  n-y\right)  -x\right)  -\left(  \left(
-y-1\right)  -x-1\right)  -1}{\left(  n-y\right)  -x}}_{\substack{=\dbinom
{n+1}{n-x-y}\\\text{(since }\left(  \left(  n-y\right)  -x\right)  -\left(
\left(  -y-1\right)  -x-1\right)  -1=n+1\\\text{and }\left(  n-y\right)
-x=n-x-y\text{)}}}\\
&  \ \ \ \ \ \ \ \ \ \ \left(
\begin{array}
[c]{c}%
\text{by (\ref{eq.binom.upper-neg}), applied to }\left(  -y-1\right)
-x-1\text{ and }\left(  n-y\right)  -x\\
\text{instead of }m\text{ and }n
\end{array}
\right) \\
&  =\left(  -1\right)  ^{n-x-y}\dbinom{n+1}{n-x-y},
\end{align*}
this yields%
\[
\left(  -1\right)  ^{n-x-y}\dbinom{n+1}{n-x-y}=\left(  -1\right)  ^{n-x-y}%
\sum_{k=0}^{n-y}\dbinom{k}{x}\dbinom{n-k}{y}.
\]
We can cancel $\left(  -1\right)  ^{n-x-y}$ from this equality (because
$\left(  -1\right)  ^{n-x-y}\neq0$). As a result, we obtain%
\begin{equation}
\dbinom{n+1}{n-x-y}=\sum_{k=0}^{n-y}\dbinom{k}{x}\dbinom{n-k}{y}.
\label{pf.prop.vandermonde.consequences.f.5}%
\end{equation}

But $0\leq n-y$ (since $n-y\in\mathbb{N}$) and $n-y\leq n$ (since
$y\in\mathbb{N}$). Also, every $k\in\left\{  n-y+1,n-y+2,\ldots,n\right\}  $
satisfies%
\begin{equation}
\dbinom{n-k}{y}=0 \label{pf.prop.vandermonde.consequences.f.2}%
\end{equation}
\footnote{\textit{Proof of (\ref{pf.prop.vandermonde.consequences.f.2}):} Let
$k\in\left\{  n-y+1,n-y+2,\ldots,n\right\}  $. Then, $k\leq n$ and $k>n-y$.
Hence, $n-k\in\mathbb{N}$ (since $k\leq n$) and $n-\underbrace{k}%
_{>n-y}<n-\left(  n-y\right)  =y$. Therefore, (\ref{eq.binom.0}) (applied to
$n-k$ and $y$ instead of $m$ and $n$) yields $\dbinom{n-k}{y}=0$. This proves
(\ref{pf.prop.vandermonde.consequences.f.2}).}. Hence,%
\begin{align}
\sum_{k=0}^{n}\dbinom{k}{x}\dbinom{n-k}{y}  &  =\sum_{k=0}^{n-y}\dbinom{k}%
{x}\dbinom{n-k}{y}+\sum_{k=n-y+1}^{n}\dbinom{k}{x}\underbrace{\dbinom{n-k}{y}%
}_{\substack{=0\\\text{(by (\ref{pf.prop.vandermonde.consequences.f.2}))}%
}}\nonumber\\
&  \ \ \ \ \ \ \ \ \ \ \left(  \text{since }0\leq n-y\leq n\right) \nonumber\\
&  =\sum_{k=0}^{n-y}\dbinom{k}{x}\dbinom{n-k}{y}+\underbrace{\sum
_{k=n-y+1}^{n}\dbinom{k}{x}0}_{=0}=\sum_{k=0}^{n-y}\dbinom{k}{x}\dbinom
{n-k}{y}\nonumber\\
&  =\dbinom{n+1}{n-x-y}\ \ \ \ \ \ \ \ \ \ \left(  \text{by
(\ref{pf.prop.vandermonde.consequences.f.5})}\right)  .
\label{pf.prop.vandermonde.consequences.f.9}%
\end{align}

Finally, $n+1\in\mathbb{N}$ and $x+y+1\in\mathbb{N}$ (since $x\in\mathbb{N}$
and $y\in\mathbb{N}$) and $\underbrace{n}_{\geq x+y}+1\geq x+y+1$. Hence,
(\ref{eq.binom.symm}) (applied to $n+1$ and $x+y+1$ instead of $m$ and $n$)
yields
\[
\dbinom{n+1}{x+y+1}=\dbinom{n+1}{\left(  n+1\right)  -\left(  x+y+1\right)
}=\dbinom{n+1}{n-x-y}%
\]
(since $\left(  n+1\right)  -\left(  x+y+1\right)  =n-x-y$). Comparing this
with (\ref{pf.prop.vandermonde.consequences.f.9}), we obtain%
\[
\dbinom{n+1}{x+y+1}=\sum_{k=0}^{n}\dbinom{k}{x}\dbinom{n-k}{y}.
\]
Thus, Proposition \ref{prop.vandermonde.consequences} \textbf{(f)} is proven
in Case 2.

We have now proven Proposition \ref{prop.vandermonde.consequences}
\textbf{(f)} in both Cases 1 and 2; thus, Proposition
\ref{prop.vandermonde.consequences} \textbf{(f)} always holds.

\textbf{(g)} Let $x\in\mathbb{Z}$ and $y\in\mathbb{N}$ and $n\in\mathbb{N}$ be
such that $x+y\geq0$ and $n\geq x$. We have $x+y\in\mathbb{N}$ (since
$x+y\geq0$). We must be in one of the following two cases:

\textit{Case 1:} We have $x+y<n$.

\textit{Case 2:} We have $x+y\geq n$.

Let us first consider Case 1. In this case, we have $x+y<n$. Thus,
$\dbinom{x+y}{n}=0$ (by (\ref{eq.binom.0}), applied to $m=x+y$). But every
$k\in\left\{  0,1,\ldots,x+y\right\}  $ satisfies $\dbinom{y}{n+k-x}%
=0$\ \ \ \ \footnote{\textit{Proof.} Let $k\in\left\{  0,1,\ldots,x+y\right\}
$. Then, $k\geq0$, so that $n+\underbrace{k}_{\geq0}-x\geq n-x>y$ (since
$n>x+y$ (since $x+y<n$)). In other words, $y<n+k-x$. Also, $n+k-x>y\geq0$, so
that $n+k-x\in\mathbb{N}$. Hence, $\dbinom{y}{n+k-x}=0$ (by (\ref{eq.binom.0}%
), applied to $y$ and $n+k-x$ instead of $m$ and $n$). Qed.}. Thus,
$\sum_{k=0}^{x+y}\dbinom{x}{k}\underbrace{\dbinom{y}{n+k-x}}_{=0}=\sum
_{k=0}^{x+y}\dbinom{x}{k}0=0$. Compared with $\dbinom{x+y}{n}=0$, this yields
$\dbinom{x+y}{n}=\sum_{k=0}^{x+y}\dbinom{x}{k}\dbinom{y}{n+k-x}$. Thus,
Proposition \ref{prop.vandermonde.consequences} \textbf{(g)} is proven in Case 1.

Let us now consider Case 2. In this case, we have $x+y\geq n$. Hence,
$\dbinom{x+y}{n}=\dbinom{x+y}{x+y-n}$ (by (\ref{eq.binom.symm}), applied to
$m=x+y$). Also, $x+y-n\in\mathbb{N}$ (since $x+y\geq n$). Therefore,
Proposition \ref{prop.vandermonde.consequences} \textbf{(a)} (applied to
$x+y-n$ instead of $n$) yields%
\[
\dbinom{x+y}{x+y-n}=\sum_{k=0}^{x+y-n}\dbinom{x}{k}\dbinom{y}{x+y-n-k}.
\]
Since $\dbinom{x+y}{n}=\dbinom{x+y}{x+y-n}$, this rewrites as
\begin{equation}
\dbinom{x+y}{n}=\sum_{k=0}^{x+y-n}\dbinom{x}{k}\dbinom{y}{x+y-n-k}.
\label{pf.prop.vandermonde.consequences.g.4}%
\end{equation}
But every $k\in\left\{  0,1,\ldots,x+y-n\right\}  $ satisfies $\dbinom
{y}{x+y-n-k}=\dbinom{y}{n+k-x}$\ \ \ \ \footnote{\textit{Proof.} Let
$k\in\left\{  0,1,\ldots,x+y-n\right\}  $. Then, $0\leq k\leq x+y-n$. Hence,
$x+y-n\geq k$, so that $x+y-n-k\in\mathbb{N}$. Also, $y\geq x+y-n-k$ (since
$y-\left(  x+y-n-k\right)  =\underbrace{n}_{\geq x}+\underbrace{k}_{\geq
0}-x\geq x+0-x=0$). Therefore, (\ref{eq.binom.symm}) (applied to $y$ and
$x+y-n-k$ instead of $m$ and $n$) yields $\dbinom{y}{x+y-n-k}=\dbinom
{y}{y-\left(  x+y-n-k\right)  }=\dbinom{y}{n+k-x}$ (since $y-\left(
x+y-n-k\right)  =n+k-x$), qed.}. Hence,
(\ref{pf.prop.vandermonde.consequences.g.4}) becomes%
\begin{equation}
\dbinom{x+y}{n}=\sum_{k=0}^{x+y-n}\dbinom{x}{k}\underbrace{\dbinom{y}%
{x+y-n-k}}_{=\dbinom{y}{n+k-x}}=\sum_{k=0}^{x+y-n}\dbinom{x}{k}\dbinom
{y}{n+k-x}. \label{pf.prop.vandermonde.consequences.g.6}%
\end{equation}

On the other hand, we have $0\leq n\leq x+y$ and thus $0\leq x+y-n\leq x+y$.
But every $k\in\mathbb{N}$ satisfying $k>x+y-n$ satisfies
\begin{equation}
\dbinom{y}{n+k-x}=0 \label{pf.prop.vandermonde.consequences.g.7}%
\end{equation}
\footnote{\textit{Proof.} Let $k\in\mathbb{N}$ be such that $k>x+y-n$. Then,
$n+\underbrace{k}_{>x+y-n}-x>n+\left(  x+y-n\right)  -x=y$. In other words,
$y<n+k-x$. Also, $n+k-x>y\geq0$, so that $n+k-x\in\mathbb{N}$. Hence,
(\ref{eq.binom.0}) (applied to $y$ and $n+k-x$ instead of $m$ and $n$) yields
$\dbinom{y}{n+k-x}=0$, qed.}. Hence,
\begin{align*}
&  \sum_{k=0}^{x+y}\dbinom{x}{k}\dbinom{y}{n+k-x}\\
&  =\sum_{k=0}^{x+y-n}\dbinom{x}{k}\dbinom{y}{n+k-x}+\sum_{k=\left(
x+y-n\right)  +1}^{x+y}\dbinom{x}{k}\underbrace{\dbinom{y}{n+k-x}%
}_{\substack{=0\\\text{(by (\ref{pf.prop.vandermonde.consequences.g.7}) (since
}k\geq\left(  x+y-n\right)  +1>x+y-n\text{))}}}\\
&  \ \ \ \ \ \ \ \ \ \ \left(  \text{since }0\leq x+y-n\leq x+y\right) \\
&  =\sum_{k=0}^{x+y-n}\dbinom{x}{k}\dbinom{y}{n+k-x}+\underbrace{\sum
_{k=\left(  x+y-n\right)  +1}^{x+y}\dbinom{x}{k}0}_{=0}=\sum_{k=0}%
^{x+y-n}\dbinom{x}{k}\dbinom{y}{n+k-x}.
\end{align*}
Compared with (\ref{pf.prop.vandermonde.consequences.g.6}), this yields%
\[
\dbinom{x+y}{n}=\sum_{k=0}^{x+y}\dbinom{x}{k}\dbinom{y}{n+k-x}.
\]
This proves Proposition \ref{prop.vandermonde.consequences} \textbf{(g)} in
Case 2.

Proposition \ref{prop.vandermonde.consequences} \textbf{(g)} is thus proven in
each of the two Cases 1 and 2. Therefore, Proposition
\ref{prop.vandermonde.consequences} \textbf{(g)} holds in full generality.
\end{proof}

\begin{remark}
The proof of Proposition \ref{prop.vandermonde.consequences} given above
illustrates a useful technique: the use of upper negation (i.e., the equality
(\ref{eq.binom.upper-neg})) to transform one equality into another. In a nutshell,

\begin{itemize}
\item we have proven Proposition \ref{prop.vandermonde.consequences}
\textbf{(d)} by applying Proposition \ref{prop.vandermonde.consequences}
\textbf{(a)} to $-y$ instead of $y$, and then rewriting the result using upper negation;

\item we have proven Proposition \ref{prop.vandermonde.consequences}
\textbf{(e)} by applying Proposition \ref{prop.vandermonde.consequences}
\textbf{(a)} to $y$, $-x-1$ and $n-x$ instead of $x$, $y$ and $n$, and then
rewriting the resulting identity using upper negation;

\item we have proven Proposition \ref{prop.vandermonde.consequences}
\textbf{(f)} by applying Proposition \ref{prop.vandermonde.consequences}
\textbf{(e)} to $-y-1$ and $n-y$ instead of $y$ and $n$, and rewriting the
resulting identity using upper negation.
\end{itemize}

Thus, by substitution and rewriting using upper negation, one single equality
(namely, Proposition \ref{prop.vandermonde.consequences} \textbf{(a)}) has
morphed into three other equalities. Note, in particular, that no negative
numbers appear in Proposition \ref{prop.vandermonde.consequences}
\textbf{(f)}, but yet we proved it by substituting negative values for $x$ and
$y$.
\end{remark}

\subsection{Further results}

\begin{exercise}
\label{exe.ps1.1.2}Let $n$ be a nonnegative integer. Prove that there exist
\textbf{nonnegative} integers $c_{i,j}$ for all $0\leq i\leq n$ and $0\leq
j\leq n$ such that%
\begin{equation}
\dbinom{XY}{n}=\sum_{i=0}^{n}\sum_{j=0}^{n}c_{i,j}\dbinom{X}{i}\dbinom{Y}{j}
\label{eq.exe.1.2.claim}%
\end{equation}
(an equality between polynomials in two variables $X$ and $Y$).
\end{exercise}

Notice that the integers $c_{i,j}$ in Exercise \ref{exe.ps1.1.2} can depend on
the $n$ (besides depending on $i$ and $j$). We just have not included the $n$
in the notation because it is fixed.

We shall now state two results that are used by Lee and Schiffler in their
celebrated proof of positivity for cluster algebras \cite{LS} (one of the
recent breakthroughs in cluster algebra theory). Specifically, our Exercise
\ref{exe.ps1.1.3} is (essentially) \cite[Lemma 5.11]{LS}, and our Proposition
\ref{prop.ps1.1.4} is (essentially) \cite[Lemma 5.12]{LS}\footnote{We say
\textquotedblleft essentially\textquotedblright\ because the $X$ in
\cite[Lemma 5.11]{LS} and in \cite[Lemma 5.12]{LS} is a variable ranging over
the nonnegative integers rather than an indeterminate. But this does not make
much of a difference (indeed, Lemma \ref{lem.polyeq} \textbf{(b)} allows us to
easily derive our Exercise \ref{exe.ps1.1.3} and Proposition
\ref{prop.ps1.1.4} from \cite[Lemma 5.11]{LS} and \cite[Lemma 5.12]{LS}, and
of course the converse implication is obvious).}.

\begin{exercise}
\label{exe.ps1.1.3}Let $a$, $b$ and $c$ be three nonnegative integers. Prove
that the polynomial $\dbinom{aX+b}{c}$ in the variable $X$ (this is a
polynomial in $X$ of degree $\leq c$) can be written as a sum $\sum_{i=0}%
^{c}d_{i}\dbinom{X}{i}$ with \textbf{nonnegative} $d_{i}$.
\end{exercise}

\begin{proposition}
\label{prop.ps1.1.4}Let $a$ and $b$ be two nonnegative integers. There exist
\textbf{nonnegative} integers $e_{0},e_{1},\ldots,e_{a+b}$ such that%
\[
\dbinom{X}{a}\dbinom{X}{b}=\sum_{i=0}^{a+b}e_{i}\dbinom{X}{i}%
\]
(an equality between polynomials in $X$).
\end{proposition}

\begin{proof}
[First proof of Proposition \ref{prop.ps1.1.4}.]For every $N\in\mathbb{N}$, we
let $\left[  N\right]  $ denote the $N$-element set $\left\{  1,2,\ldots
,N\right\}  $.

For every set $S$, we let an $S$\textit{-junction} mean a pair $\left(
A,B\right)  $, where $A$ is an $a$-element subset of $S$ and where $B$ is a
$b$-element subset of $S$ such that $A\cup B=S$. (We do not mention $a$ and
$b$ in our notation, because $a$ and $b$ are fixed.)

For example, if $a=2$ and $b=3$, then $\left(  \left\{  1,4\right\}  ,\left\{
2,3,4\right\}  \right)  $ is a $\left[  4\right]  $-junction, and $\left(
\left\{  2,4\right\}  ,\left\{  1,4,6\right\}  \right)  $ is a $\left\{
1,2,4,6\right\}  $-junction, but $\left(  \left\{  1,3\right\}  ,\left\{
2,3,5\right\}  \right)  $ is not a $\left[  5\right]  $-junction (since
$\left\{  1,3\right\}  \cup\left\{  2,3,5\right\}  \neq\left[  5\right]  $).

For every $i\in\mathbb{N}$, we let $e_{i}$ be the number of all $\left[
i\right]  $-junctions. Then, if $S$ is any $i$-element set, then%
\begin{equation}
e_{i}\text{ is the number of all }S\text{-junctions} \label{pf.prop.ps1.1.4.1}%
\end{equation}
\footnote{\textit{Proof of (\ref{pf.prop.ps1.1.4.1}):} Let $S$ be any
$i$-element set. We know that $e_{i}$ is the number of all $\left[  i\right]
$-junctions. We want to prove that $e_{i}$ is the number of all $S$-junctions.
Roughly speaking, this is obvious, because we can \textquotedblleft relabel
the elements of $S$ as $1,2,\ldots,i$\textquotedblright\ (since $S$ is an
$i$-element set), and then the $S$-junctions become precisely the $\left[
i\right]  $-junctions.
\par
Here is a formal way to make this argument: The sets $\left[  i\right]  $ and
$S$ have the same number of elements (indeed, both are $i$-element sets).
Hence, there exists a bijection $\phi:S\rightarrow\left[  i\right]  $. Fix
such a $\phi$. Now, the $S$-junctions are in a 1-to-1 correspondence with the
$\left[  i\right]  $-junctions (namely, to every $S$-junction $\left(
A,B\right)  $ corresponds the $\left[  i\right]  $-junction $\left(
\phi\left(  A\right)  ,\phi\left(  B\right)  \right)  $, and conversely, to
every $\left[  i\right]  $-junction $\left(  A^{\prime},B^{\prime}\right)  $
corresponds the $S$-junction $\left(  \phi^{-1}\left(  A^{\prime}\right)
,\phi^{-1}\left(  B^{\prime}\right)  \right)  $). Hence, the number of all
$S$-junctions equals the number of $\left[  i\right]  $-junctions. Since the
latter number is $e_{i}$, this shows that the former number is also $e_{i}$.
This proves (\ref{pf.prop.ps1.1.4.1}).}.

Now, let us show that%
\begin{equation}
\dbinom{x}{a}\dbinom{x}{b}=\sum_{i=0}^{a+b}e_{i}\dbinom{x}{i}
\label{pf.prop.ps1.1.4.xclaim}%
\end{equation}
for every $x\in\mathbb{N}$.

[\textit{Proof of (\ref{pf.prop.ps1.1.4.xclaim}):} Let $x\in\mathbb{N}$. How
many ways are there to choose a pair $\left(  A,B\right)  $ consisting of an
$a$-element subset $A$ of $\left[  x\right]  $ and a $b$-element subset $B$ of
$\left[  x\right]  $ ?

Let us give two different answers to this question. The first answer is the
straightforward one: To choose a pair $\left(  A,B\right)  $ consisting of an
$a$-element subset $A$ of $\left[  x\right]  $ and a $b$-element subset $B$ of
$\left[  x\right]  $, we need to choose an $a$-element subset $A$ of $\left[
x\right]  $ and a $b$-element subset $B$ of $\left[  x\right]  $. There are
$\dbinom{x}{a}\dbinom{x}{b}$ total ways to do this (since there are
$\dbinom{x}{a}$ choices for $A$\ \ \ \ \footnote{This follows from
(\ref{eq.binom.subsets}).}, and $\dbinom{x}{b}$ choices for $B$%
\ \ \ \ \footnote{Again, this follows from (\ref{eq.binom.subsets}).}, and
these choices are independent). In other words, the number of all pairs
$\left(  A,B\right)  $ consisting of an $a$-element subset $A$ of $\left[
x\right]  $ and a $b$-element subset $B$ of $\left[  x\right]  $ equals
$\dbinom{x}{a}\dbinom{x}{b}$.

On the other hand, here is a more imaginative procedure to choose a pair
$\left(  A,B\right)  $ consisting of an $a$-element subset $A$ of $\left[
x\right]  $ and a $b$-element subset $B$ of $\left[  x\right]  $:

\begin{enumerate}
\item We choose how many elements the union $A\cup B$ will have. In other
words, we choose an $i\in\mathbb{N}$ that will satisfy $\left\vert A\cup
B\right\vert =i$. This $i$ must be an integer between $0$ and $a+b$
(inclusive)\footnote{\textit{Proof.} Clearly, $i$ cannot be smaller than $0$.
But $i$ also cannot be larger than $a+b$ (since $i$ will have to satisfy
$i=\left\vert A\cup B\right\vert \leq\underbrace{\left\vert A\right\vert
}_{=a}+\underbrace{\left\vert B\right\vert }_{=b}=a+b$). Thus, $i$ must be an
integer between $0$ and $a+b$ (inclusive).}.

\item We choose a subset $S$ of $\left[  x\right]  $, which will serve as the
union $A\cup B$. This subset $S$ must be an $i$-element subset of $\left[
x\right]  $ (because we will have $\left\vert \underbrace{S}_{=A\cup
B}\right\vert =\left\vert A\cup B\right\vert =i$). Thus, there are $\dbinom
{x}{i}$ ways to choose it (since we need to choose an $i$-element subset of
$\left[  x\right]  $).

\item Now, it remains to choose the pair $\left(  A,B\right)  $ itself. This
pair must be a pair of subsets of $\left[  x\right]  $ satisfying $\left\vert
A\right\vert =a$, $\left\vert B\right\vert =b$, $A\cup B=S$ and $\left\vert
A\cup B\right\vert =i$. We can forget about the $\left\vert A\cup B\right\vert
=i$ condition, since it automatically follows from $A\cup B=S$ (because
$\left\vert S\right\vert =i$). So we need to choose a pair $\left(
A,B\right)  $ of subsets of $\left[  x\right]  $ satisfying $\left\vert
A\right\vert =a$, $\left\vert B\right\vert =b$ and $A\cup B=S$. In other
words, we need to choose a pair $\left(  A,B\right)  $ of subsets of $S$
satisfying $\left\vert A\right\vert =a$, $\left\vert B\right\vert =b$ and
$A\cup B=S$\ \ \ \ \footnote{Here, we have replaced \textquotedblleft subsets
of $\left[  x\right]  $\textquotedblright\ by \textquotedblleft subsets of
$S$\textquotedblright, because the condition $A\cup B=S$ forces $A$ and $B$ to
be subsets of $S$.}. In other words, we need to choose an $S$-junction (since
this is how an $S$-junction was defined). This can be done in exactly $e_{i}$
ways (according to (\ref{pf.prop.ps1.1.4.1})).
\end{enumerate}

Thus, in total, there are $\sum_{i=0}^{a+b}\dbinom{x}{i}e_{i}$ ways to perform
this procedure. Hence, the total number of all pairs $\left(  A,B\right)  $
consisting of an $a$-element subset $A$ of $\left[  x\right]  $ and a
$b$-element subset $B$ of $\left[  x\right]  $ equals $\sum_{i=0}^{a+b}%
\dbinom{x}{i}e_{i}$. But earlier, we have shown that this number is
$\dbinom{x}{a}\dbinom{x}{b}$. Comparing these two results, we conclude that
$\dbinom{x}{a}\dbinom{x}{b}=\sum_{i=0}^{a+b}\dbinom{x}{i}e_{i}=\sum
_{i=0}^{a+b}e_{i}\dbinom{x}{i}$. Thus, (\ref{pf.prop.ps1.1.4.xclaim}) is proven.]

Now, we define two polynomials $P$ and $Q$ in the indeterminate $X$ with
rational coefficients by setting%
\[
P=\dbinom{X}{a}\dbinom{X}{b};\ \ \ \ \ \ \ \ \ \ Q=\sum_{i=0}^{a+b}%
e_{i}\dbinom{X}{i}.
\]
The equality (\ref{pf.prop.ps1.1.4.xclaim}) (which we have proven) states that
$P\left(  x\right)  =Q\left(  x\right)  $ for all $x\in\mathbb{N}$. Thus,
Lemma \ref{lem.polyeq} \textbf{(b)} yields that $P=Q$. Recalling how $P$ and
$Q$ are defined, we see that this rewrites as $\dbinom{X}{a}\dbinom{X}{b}%
=\sum_{i=0}^{a+b}e_{i}\dbinom{X}{i}$. This proves Proposition
\ref{prop.ps1.1.4}.
\end{proof}

Our second proof of Proposition \ref{prop.ps1.1.4} is algebraic, and is based
on a suggestion of math.stackexchange user tcamps in a comment on
\href{http://math.stackexchange.com/questions/1342384}{question \#1342384}. It
proceeds by way of the following, more explicit result:

\begin{proposition}
\label{prop.binom.XaXb.X}Let $a$ and $b$ be two nonnegative integers. Then,%
\[
\dbinom{X}{a}\dbinom{X}{b}=\sum_{i=a}^{a+b}\dbinom{i}{a}\dbinom{a}%
{a+b-i}\dbinom{X}{i}.
\]

\end{proposition}

Let us also state the analogue of this proposition in which the indeterminate
$X$ is replaced by a rational number $m$:

\begin{proposition}
\label{prop.binom.XaXb.rat}Let $a$ and $b$ be two nonnegative integers. Let
$m\in\mathbb{Q}$. Then,%
\[
\dbinom{m}{a}\dbinom{m}{b}=\sum_{i=a}^{a+b}\dbinom{i}{a}\dbinom{a}%
{a+b-i}\dbinom{m}{i}.
\]

\end{proposition}

\begin{proof}
[Proof of Proposition \ref{prop.binom.XaXb.rat}.]Theorem
\ref{thm.vandermonde.rat} (applied to $b$, $m-a$ and $a$ instead of $n$, $x$
and $y$) yields%
\[
\dbinom{\left(  m-a\right)  +a}{b}=\sum_{k=0}^{b}\dbinom{m-a}{k}\dbinom
{a}{b-k}.
\]
Since $\left(  m-a\right)  +a=m$, this rewrites as%
\begin{align*}
\dbinom{m}{b}  &  =\sum_{k=0}^{b}\dbinom{m-a}{k}\dbinom{a}{b-k}=\sum
_{i=a}^{a+b}\dbinom{m-a}{i-a}\underbrace{\dbinom{a}{b-\left(  i-a\right)  }%
}_{\substack{=\dbinom{a}{a+b-i}\\\text{(since }b-\left(  i-a\right)
=a+b-i\text{)}}}\\
&  \ \ \ \ \ \ \ \ \ \ \left(  \text{here, we substituted }i-a\text{ for
}k\text{ in the sum}\right) \\
&  =\sum_{i=a}^{a+b}\dbinom{m-a}{i-a}\dbinom{a}{a+b-i}.
\end{align*}
Multiplying both sides of this identity with $\dbinom{m}{a}$, we obtain%
\begin{align*}
\dbinom{m}{a}\dbinom{m}{b}  &  =\dbinom{m}{a}\sum_{i=a}^{a+b}\dbinom{m-a}%
{i-a}\dbinom{a}{a+b-i}=\sum_{i=a}^{a+b}\underbrace{\dbinom{m}{a}\dbinom
{m-a}{i-a}}_{\substack{=\dbinom{m}{i}\dbinom{i}{a}\\\text{(by Proposition
\ref{prop.binom.trinom-rev})}}}\dbinom{a}{a+b-i}\\
&  =\sum_{i=a}^{a+b}\dbinom{m}{i}\dbinom{i}{a}\dbinom{a}{a+b-i}=\sum
_{i=a}^{a+b}\dbinom{i}{a}\dbinom{a}{a+b-i}\dbinom{m}{i}.
\end{align*}
This proves Proposition \ref{prop.binom.XaXb.rat}.
\end{proof}

\begin{proof}
[Proof of Proposition \ref{prop.binom.XaXb.X}.]To obtain a proof of
Proposition \ref{prop.binom.XaXb.X}, replace every appearance of
\textquotedblleft$m$\textquotedblright\ by \textquotedblleft$X$%
\textquotedblright\ in the above proof of Proposition
\ref{prop.binom.XaXb.rat}. (Of course, this requires knowing that Theorem
\ref{thm.vandermonde.rat} holds when $x$ and $y$ are polynomials rather than
numbers. But this is true, because both proofs that we gave for Theorem
\ref{thm.vandermonde.rat} still apply in this case.)
\end{proof}

\begin{proof}
[Second proof of Proposition \ref{prop.ps1.1.4}.]Let us define $a+b+1$
nonnegative integers \newline$e_{0},e_{1},\ldots,e_{a+b}$ by%
\begin{equation}
e_{i}=%
\begin{cases}
\dbinom{i}{a}\dbinom{a}{a+b-i}, & \text{if }i\geq a;\\
0, & \text{otherwise}%
\end{cases}
\ \ \ \ \ \ \ \ \ \ \text{for all }i\in\left\{  0,1,\ldots,a+b\right\}  .
\label{pf.prop.ps1.1.4.sol2.ei}%
\end{equation}
Then,%
\begin{align*}
\sum_{i=0}^{a+b}e_{i}\dbinom{X}{i}  &  =\sum_{i=a}^{a+b}\dbinom{i}{a}%
\dbinom{a}{a+b-i}\dbinom{X}{i}\ \ \ \ \ \ \ \ \ \ \left(  \text{by our
definition of }e_{0},e_{1},\ldots,e_{a+b}\right) \\
&  =\dbinom{X}{a}\dbinom{X}{b}\ \ \ \ \ \ \ \ \ \ \left(  \text{by Proposition
\ref{prop.binom.XaXb.X}}\right)  .
\end{align*}
Thus, Proposition \ref{prop.ps1.1.4} is proven again.
\end{proof}

\begin{remark}
Comparing our two proofs of Proposition \ref{prop.ps1.1.4}, it is natural to
suspect that the $e_{0},e_{1},\ldots,e_{a+b}$ defined in the First proof are
identical with the $e_{0},e_{1},\ldots,e_{a+b}$ defined in the Second proof.
This actually follows from general principles (namely, from the word
\textquotedblleft unique\textquotedblright\ in Proposition
\ref{prop.hartshorne} \textbf{(a)}), but there is also a simple combinatorial
reason. Namely, let $i\in\left\{  0,1,\ldots,a+b\right\}  $. We shall show
that the $e_{i}$ defined in the First proof equals the $e_{i}$ defined in the
Second proof.

The $e_{i}$ defined in the First proof is the number of all $\left[  i\right]
$-junctions. An $\left[  i\right]  $-junction is a pair $\left(  A,B\right)
$, where $A$ is an $a$-element subset of $\left[  i\right]  $ and where $B$ is
a $b$-element subset of $\left[  i\right]  $ such that $A\cup B=\left[
i\right]  $. Here is a way to construct an $\left[  i\right]  $-junction:

\begin{itemize}
\item First, we pick the set $A$. There are $\dbinom{i}{a}$ ways to do this,
since $A$ has to be an $a$-element subset of the $i$-element set $\left[
i\right]  $.

\item Then, we pick the set $B$. This has to be a $b$-element subset of the
$i$-element set $\left[  i\right]  $ satisfying $A\cup B=\left[  i\right]  $.
The equality $A\cup B=\left[  i\right]  $ means that $B$ has to contain the
$i-a$ elements of $\left[  i\right]  \setminus A$; but the remaining
$b-\left(  i-a\right)  =a+b-i$ elements of $B$ can be chosen arbitrarily among
the $a$ elements of $A$. Thus, there are $\dbinom{a}{a+b-i}$ ways to choose
$B$ (since we have to choose $a+b-i$ elements of $B$ among the $a$ elements of
$A$).
\end{itemize}

Thus, the number of all $\left[  i\right]  $-junctions is $\dbinom{i}%
{a}\dbinom{a}{a+b-i}$. This can be rewritten in the form $%
\begin{cases}
\dbinom{i}{a}\dbinom{a}{a+b-i}, & \text{if }i\geq a;\\
0, & \text{otherwise}%
\end{cases}
$ (because if $i<a$, then $\dbinom{i}{a}=0$ and thus $\dbinom{i}{a}\dbinom
{a}{a+b-i}=0$). Thus, we have shown that the number of all $\left[  i\right]
$-junctions is $%
\begin{cases}
\dbinom{i}{a}\dbinom{a}{a+b-i}, & \text{if }i\geq a;\\
0, & \text{otherwise}%
\end{cases}
$. In other words, the $e_{i}$ defined in the First proof equals the $e_{i}$
defined in the Second proof.
\end{remark}

Here is an assortment of other identities that involve binomial coefficients:

\begin{proposition}
\label{prop.binom.bin-id}\textbf{(a)} Every $x\in\mathbb{Z}$, $y\in\mathbb{Z}$
and $n\in\mathbb{N}$ satisfy $\left(  x+y\right)  ^{n}=\sum_{k=0}^{n}%
\dbinom{n}{k}x^{k}y^{n-k}$.

\textbf{(b)} Every $n\in\mathbb{N}$ satisfies $\sum_{k=0}^{n}\dbinom{n}%
{k}=2^{n}$.

\textbf{(c)} Every $n\in\mathbb{N}$ satisfies $\sum_{k=0}^{n}\left(
-1\right)  ^{k}\dbinom{n}{k}=%
\begin{cases}
1, & \text{if }n=0;\\
0, & \text{if }n\neq0
\end{cases}
$.

\textbf{(d)} Every $n\in\mathbb{Z}$, $i\in\mathbb{N}$ and $a\in\mathbb{N}$
satisfying $i\geq a$ satisfy $\dbinom{n}{i}\dbinom{i}{a}=\dbinom{n}{a}%
\dbinom{n-a}{i-a}$.

\textbf{(e)} Every $n\in\mathbb{N}$ and $m\in\mathbb{Z}$ satisfy $\sum
_{i=0}^{n}\dbinom{n}{i}\dbinom{m+i}{n}=\sum_{i=0}^{n}\dbinom{n}{i}\dbinom
{m}{i}2^{i}$.

\textbf{(f)} Every $a\in\mathbb{N}$, $b\in\mathbb{N}$ and $x\in\mathbb{Z}$
satisfy $\sum_{i=0}^{b}\dbinom{a}{i}\dbinom{b}{i}\dbinom{x+i}{a+b}=\dbinom
{x}{a}\dbinom{x}{b}$.

\textbf{(g)} Every $a\in\mathbb{N}$, $b\in\mathbb{N}$ and $x\in\mathbb{Z}$
satisfy $\sum_{i=0}^{b}\dbinom{a}{i}\dbinom{b}{i}\dbinom{a+b+x-i}{a+b}%
=\dbinom{a+x}{a}\dbinom{b+x}{b}$.
\end{proposition}

(I have learnt parts \textbf{(e)} and \textbf{(f)} of Proposition
\ref{prop.binom.bin-id} from
\href{http://artofproblemsolving.com/community/c6h626709p4577721}{AoPS}, but
they are fairly classical results. Part \textbf{(e)} is equivalent to a claim
in \cite[Chapter I, Exercise 21]{Comtet74}. Part \textbf{(f)} is \cite[\S 1.4,
(10)]{Riorda68}. Part \textbf{(g)} is a restatement of \cite[(6.93)]{Gould-I}.)

\begin{proof}
[Proof of Proposition \ref{prop.binom.bin-id}.]\textbf{(a)} Proposition
\ref{prop.binom.bin-id} \textbf{(a)} is clearly a particular case of
Proposition \ref{prop.binom.binomial}.

\textbf{(b)} Let $n\in\mathbb{N}$. Applying Proposition
\ref{prop.binom.bin-id} \textbf{(a)} to $x=1$ and $y=1$, we obtain%
\[
\left(  1+1\right)  ^{n}=\sum_{k=0}^{n}\dbinom{n}{k}\underbrace{1^{k}}%
_{=1}\underbrace{1^{n-k}}_{=1}=\sum_{k=0}^{n}\dbinom{n}{k},
\]
thus%
\[
\sum_{k=0}^{n}\dbinom{n}{k}=\left(  \underbrace{1+1}_{=2}\right)  ^{n}=2^{n}.
\]
This proves Proposition \ref{prop.binom.bin-id} \textbf{(b)}.

\textbf{(c)} Let $n\in\mathbb{N}$. Applying Proposition
\ref{prop.binom.bin-id} \textbf{(a)} to $x=-1$ and $y=1$, we obtain
\[
\left(  -1+1\right)  ^{n}=\sum_{k=0}^{n}\dbinom{n}{k}\left(  -1\right)
^{k}\underbrace{1^{n-k}}_{=1}=\sum_{k=0}^{n}\dbinom{n}{k}\left(  -1\right)
^{k}=\sum_{k=0}^{n}\left(  -1\right)  ^{k}\dbinom{n}{k},
\]
thus%
\[
\sum_{k=0}^{n}\left(  -1\right)  ^{k}\dbinom{n}{k}=\left(  \underbrace{-1+1}%
_{=0}\right)  ^{n}=0^{n}=%
\begin{cases}
1, & \text{if }n=0;\\
0, & \text{if }n\neq0
\end{cases}
.
\]
This proves Proposition \ref{prop.binom.bin-id} \textbf{(c)}.

\textbf{(d)} Let $n\in\mathbb{Z}$, $i\in\mathbb{N}$ and $a\in\mathbb{N}$ be
such that $i\geq a$. Proposition \ref{prop.binom.trinom-rev} (applied to $n$
instead of $m$) yields $\dbinom{n}{i}\dbinom{i}{a}=\dbinom{n}{a}\dbinom
{n-a}{i-a}$. This proves Proposition \ref{prop.binom.bin-id} \textbf{(d)}.

\textbf{(e)} Let $n\in\mathbb{N}$ and $m\in\mathbb{Z}$. Clearly, every
$p\in\mathbb{N}$ satisfies%
\begin{align}
\sum_{i=0}^{p}\dbinom{p}{i}  &  =\sum_{k=0}^{p}\dbinom{p}{k}%
\ \ \ \ \ \ \ \ \ \ \left(  \text{here, we renamed the summation index
}i\text{ as }k\right) \nonumber\\
&  =2^{p} \label{pf.prop.binom.bin-id.e.4}%
\end{align}
(by Proposition \ref{prop.binom.bin-id} \textbf{(b)}, applied to $p$ instead
of $n$).

Now, let $i\in\left\{  0,1,\ldots,n\right\}  $. Applying Proposition
\ref{prop.vandermonde.consequences} \textbf{(a)} to $x=i$ and $y=m$, we obtain%
\begin{align}
&  \dbinom{i+m}{n}\nonumber\\
&  =\sum_{k=0}^{n}\dbinom{i}{k}\dbinom{m}{n-k}\nonumber\\
&  =\sum_{k=0}^{i}\dbinom{i}{k}\dbinom{m}{n-k}+\sum_{k=i+1}^{n}%
\underbrace{\dbinom{i}{k}}_{\substack{=0\\\text{(by (\ref{eq.binom.0}),
applied to }i\text{ and }k\\\text{instead of }m\text{ and }n\text{ (since
}i<k\\\text{(because }k\geq i+1>i\text{)))}}}\dbinom{m}{n-k}%
\ \ \ \ \ \ \ \ \ \ \left(  \text{since }0\leq i\leq n\right) \nonumber\\
&  =\sum_{k=0}^{i}\dbinom{i}{k}\dbinom{m}{n-k}+\underbrace{\sum_{k=i+1}%
^{n}0\dbinom{m}{n-k}}_{=0}=\sum_{k=0}^{i}\dbinom{i}{k}\dbinom{m}{n-k}.
\label{pf.prop.binom.bin-id.e.1}%
\end{align}

Now, let us forget that we fixed $i$. We thus have proven
(\ref{pf.prop.binom.bin-id.e.1}) for every $i\in\left\{  0,1,\ldots,n\right\}
$. Now,%
\begin{align*}
&  \sum_{i=0}^{n}\dbinom{n}{i}\underbrace{\dbinom{m+i}{n}}_{\substack{=\dbinom
{i+m}{n}=\sum_{k=0}^{i}\dbinom{i}{k}\dbinom{m}{n-k}\\\text{(by
(\ref{pf.prop.binom.bin-id.e.1}))}}}\\
&  =\sum_{i=0}^{n}\dbinom{n}{i}\left(  \sum_{k=0}^{i}\dbinom{i}{k}\dbinom
{m}{n-k}\right)  =\underbrace{\sum_{i=0}^{n}\sum_{k=0}^{i}}_{=\sum_{k=0}%
^{n}\sum_{i=k}^{n}}\dbinom{n}{i}\dbinom{i}{k}\dbinom{m}{n-k}\\
&  =\sum_{k=0}^{n}\sum_{i=k}^{n}\underbrace{\dbinom{n}{i}\dbinom{i}{k}%
}_{\substack{=\dbinom{n}{k}\dbinom{n-k}{i-k}\\\text{(by Proposition
\ref{prop.binom.bin-id} \textbf{(d)},}\\\text{applied to }a=k\text{ (since
}i\geq k\text{))}}}\dbinom{m}{n-k}\\
&  =\sum_{k=0}^{n}\sum_{i=k}^{n}\dbinom{n}{k}\dbinom{n-k}{i-k}\dbinom{m}%
{n-k}=\sum_{k=0}^{n}\dbinom{n}{k}\dbinom{m}{n-k}\sum_{i=k}^{n}\dbinom
{n-k}{i-k}\\
&  =\sum_{k=0}^{n}\underbrace{\dbinom{n}{k}}_{\substack{=\dbinom{n}%
{n-k}\\\text{(by (\ref{eq.binom.symm}), applied to }n\text{ and }%
k\\\text{instead of }m\text{ and }n\text{ (since }n\geq k\text{))}}}\dbinom
{m}{n-k}\sum_{i=0}^{n-k}\dbinom{n-k}{i}\\
&  \ \ \ \ \ \ \ \ \ \ \left(  \text{here, we have substituted }i\text{ for
}i-k\text{ in the second sum}\right) \\
&  =\sum_{k=0}^{n}\dbinom{n}{n-k}\dbinom{m}{n-k}\sum_{i=0}^{n-k}\dbinom
{n-k}{i}=\sum_{k=0}^{n}\dbinom{n}{k}\dbinom{m}{k}\underbrace{\sum_{i=0}%
^{k}\dbinom{k}{i}}_{\substack{=2^{k}\\\text{(by
(\ref{pf.prop.binom.bin-id.e.4}), applied to }p=k\text{)}}}\\
&  \ \ \ \ \ \ \ \ \ \ \left(  \text{here, we have substituted }k\text{ for
}n-k\text{ in the first sum}\right) \\
&  =\sum_{k=0}^{n}\dbinom{n}{k}\dbinom{m}{k}2^{k}=\sum_{i=0}^{n}\dbinom{n}%
{i}\dbinom{m}{i}2^{i}%
\end{align*}
(here, we have renamed the summation index $k$ as $i$). This proves
Proposition \ref{prop.binom.bin-id} \textbf{(e)}.

\textbf{(f)} Let $a\in\mathbb{N}$, $b\in\mathbb{N}$ and $x\in\mathbb{Z}$.
Proposition \ref{prop.binom.XaXb.rat} (applied to $m=x$) yields%
\begin{equation}
\dbinom{x}{a}\dbinom{x}{b}=\sum_{i=a}^{a+b}\dbinom{i}{a}\dbinom{a}%
{a+b-i}\dbinom{x}{i}. \label{pf.prop.binom.bin-id.f.1}%
\end{equation}

Clearly,%
\begin{align}
&  \sum_{i=0}^{a+b}\dbinom{a}{i}\dbinom{b}{i}\dbinom{x+i}{a+b}\nonumber\\
&  =\sum_{i=0}^{b}\dbinom{a}{i}\dbinom{b}{i}\dbinom{x+i}{a+b}+\sum
_{i=b+1}^{a+b}\dbinom{a}{i}\underbrace{\dbinom{b}{i}}_{\substack{=0\\\text{(by
(\ref{eq.binom.0}), applied to }b\text{ and }i\\\text{instead of }m\text{ and
}n\text{ (since }b<i\\\text{(because }i\geq b+1>b\text{)))}}}\dbinom{x+i}%
{a+b}\nonumber\\
&  \ \ \ \ \ \ \ \ \ \ \left(  \text{since }0\leq b\leq a+b\right) \nonumber\\
&  =\sum_{i=0}^{b}\dbinom{a}{i}\dbinom{b}{i}\dbinom{x+i}{a+b}+\underbrace{\sum
_{i=b+1}^{a+b}\dbinom{a}{i}0\dbinom{x+i}{a+b}}_{=0}\nonumber\\
&  =\sum_{i=0}^{b}\dbinom{a}{i}\dbinom{b}{i}\dbinom{x+i}{a+b}.
\label{pf.prop.binom.bin-id.f.2}%
\end{align}

For every $i\in\left\{  0,1,\ldots,b\right\}  $, we have%
\begin{equation}
\dbinom{x+i}{a+b}=\sum_{k=0}^{a+b}\dbinom{x}{k}\dbinom{i}{a+b-k}.
\label{pf.prop.binom.bin-id.f.2b}%
\end{equation}
(This follows from Theorem \ref{thm.vandermonde.rat} (applied to $a+b$ and $i$
instead of $n$ and $y$).) Hence,%
\begin{align*}
&  \sum_{i=0}^{a+b}\dbinom{a}{i}\dbinom{b}{i}\underbrace{\dbinom{x+i}{a+b}%
}_{\substack{=\sum_{k=0}^{a+b}\dbinom{x}{k}\dbinom{i}{a+b-k}\\\text{(by
(\ref{pf.prop.binom.bin-id.f.2b}))}}}\\
&  =\sum_{i=0}^{a+b}\dbinom{a}{i}\dbinom{b}{i}\sum_{k=0}^{a+b}\dbinom{x}%
{k}\dbinom{i}{a+b-k}\\
&  =\underbrace{\sum_{i=0}^{a+b}\sum_{k=0}^{a+b}}_{=\sum_{k=0}^{a+b}\sum
_{i=0}^{a+b}}\dbinom{a}{i}\underbrace{\dbinom{b}{i}\dbinom{x}{k}\dbinom
{i}{a+b-k}}_{=\dbinom{i}{a+b-k}\dbinom{b}{i}\dbinom{x}{k}}\\
&  =\sum_{k=0}^{a+b}\sum_{i=0}^{a+b}\dbinom{a}{i}\dbinom{i}{a+b-k}\dbinom
{b}{i}\dbinom{x}{k}=\sum_{k=0}^{a+b}\dbinom{x}{k}\sum_{i=0}^{a+b}\dbinom{a}%
{i}\dbinom{i}{a+b-k}\dbinom{b}{i}.
\end{align*}
Compared with (\ref{pf.prop.binom.bin-id.f.2}), this yields%
\begin{align}
&  \sum_{i=0}^{b}\dbinom{a}{i}\dbinom{b}{i}\dbinom{x+i}{a+b}\nonumber\\
&  =\sum_{k=0}^{a+b}\dbinom{x}{k}\sum_{i=0}^{a+b}\dbinom{a}{i}\dbinom
{i}{a+b-k}\dbinom{b}{i}. \label{pf.prop.binom.bin-id.f.3}%
\end{align}

However, for every $k\in\left\{  0,1,\ldots,a+b\right\}  $, we have%
\begin{equation}
\sum_{i=0}^{a+b}\dbinom{a}{i}\dbinom{i}{a+b-k}\dbinom{b}{i}=\dbinom{a}%
{a+b-k}\sum_{j=0}^{k}\dbinom{k-b}{k-j}\dbinom{b}{a+b-j}.
\label{pf.prop.binom.bin-id.f.4}%
\end{equation}

[\textit{Proof of (\ref{pf.prop.binom.bin-id.f.4}):} Let $k\in\left\{
0,1,\ldots,a+b\right\}  $. Then, $a+b-k\in\left\{  0,1,\ldots,a+b\right\}  $,
so that $0\leq a+b-k\leq a+b$. Now,%
\begin{align*}
&  \sum_{i=0}^{a+b}\dbinom{a}{i}\dbinom{i}{a+b-k}\dbinom{b}{i}\\
&  =\sum_{i=0}^{\left(  a+b-k\right)  -1}\dbinom{a}{i}\underbrace{\dbinom
{i}{a+b-k}}_{\substack{=0\\\text{(by (\ref{eq.binom.0}), applied to }i\text{
and}\\a+b-k\text{ instead of }m\text{ and }n\\\text{(since }i\leq\left(
a+b-k\right)  -1<a+b-k\text{))}}}\dbinom{b}{i}+\sum_{i=a+b-k}^{a+b}\dbinom
{a}{i}\dbinom{i}{a+b-k}\dbinom{b}{i}\\
&  \ \ \ \ \ \ \ \ \ \ \left(  \text{since }0\leq a+b-k\leq a+b\right) \\
&  =\underbrace{\sum_{i=0}^{\left(  a+b-k\right)  -1}\dbinom{a}{i}0\dbinom
{b}{i}}_{=0}+\sum_{i=a+b-k}^{a+b}\dbinom{a}{i}\dbinom{i}{a+b-k}\dbinom{b}{i}\\
&  =\sum_{i=a+b-k}^{a+b}\underbrace{\dbinom{a}{i}\dbinom{i}{a+b-k}%
}_{\substack{=\dbinom{a}{a+b-k}\dbinom{a-\left(  a+b-k\right)  }{i-\left(
a+b-k\right)  }\\\text{(by Proposition \ref{prop.binom.bin-id} \textbf{(d)},
applied to}\\a\text{ and }a+b-k\text{ instead of }n\text{ and }a\text{ (since
}i\geq a+b-k\text{))}}}\dbinom{b}{i}\\
&  =\sum_{i=a+b-k}^{a+b}\dbinom{a}{a+b-k}\dbinom{a-\left(  a+b-k\right)
}{i-\left(  a+b-k\right)  }\dbinom{b}{i}\\
&  =\sum_{j=0}^{k}\dbinom{a}{a+b-k}\underbrace{\dbinom{a-\left(  a+b-k\right)
}{\left(  a+b-j\right)  -\left(  a+b-k\right)  }}_{\substack{=\dbinom
{k-b}{k-j}\\\text{(since }a-\left(  a+b-k\right)  =k-b\text{ and}\\\left(
a+b-j\right)  -\left(  a+b-k\right)  =k-j\text{)}}}\dbinom{b}{a+b-j}\\
&  \ \ \ \ \ \ \ \ \ \ \left(  \text{here, we have substituted }a+b-j\text{
for }i\text{ in the sum}\right) \\
&  =\sum_{j=0}^{k}\dbinom{a}{a+b-k}\dbinom{k-b}{k-j}\dbinom{b}{a+b-j}%
=\dbinom{a}{a+b-k}\sum_{j=0}^{k}\dbinom{k-b}{k-j}\dbinom{b}{a+b-j},
\end{align*}
and this proves (\ref{pf.prop.binom.bin-id.f.4}).]

Now, (\ref{pf.prop.binom.bin-id.f.3}) becomes%
\begin{align}
&  \sum_{i=0}^{b}\dbinom{a}{i}\dbinom{b}{i}\dbinom{x+i}{a+b}\nonumber\\
&  =\sum_{k=0}^{a+b}\dbinom{x}{k}\underbrace{\sum_{i=0}^{a+b}\dbinom{a}%
{i}\dbinom{i}{a+b-k}\dbinom{b}{i}}_{\substack{=\dbinom{a}{a+b-k}\sum_{j=0}%
^{k}\dbinom{k-b}{k-j}\dbinom{b}{a+b-j}\\\text{(by
(\ref{pf.prop.binom.bin-id.f.4}))}}}\nonumber\\
&  =\sum_{k=0}^{a+b}\dbinom{x}{k}\dbinom{a}{a+b-k}\sum_{j=0}^{k}\dbinom
{k-b}{k-j}\dbinom{b}{a+b-j}\nonumber\\
&  =\sum_{i=0}^{a+b}\dbinom{x}{i}\dbinom{a}{a+b-i}\sum_{j=0}^{i}\dbinom
{i-b}{i-j}\dbinom{b}{a+b-j} \label{pf.prop.binom.bin-id.f.7}%
\end{align}
(here, we renamed the summation index $k$ as $i$ in the first sum).

Furthermore, every $i\in\left\{  0,1,\ldots,a+b\right\}  $ satisfies
\begin{equation}
\sum_{j=0}^{i}\dbinom{i-b}{i-j}\dbinom{b}{a+b-j}=\dbinom{i}{a}.
\label{pf.prop.binom.bin-id.f.8}%
\end{equation}

[\textit{Proof of (\ref{pf.prop.binom.bin-id.f.8}):} Let $i\in\left\{
0,1,\ldots,a+b\right\}  $. Thus, $0\leq i\leq a+b$. We have
\begin{align}
\sum_{j=0}^{i}\dbinom{i-b}{i-j}\dbinom{b}{a+b-j}  &  =\sum_{k=0}%
^{i}\underbrace{\dbinom{i-b}{i-\left(  i-k\right)  }}_{\substack{=\dbinom
{i-b}{k}\\\text{(since }i-\left(  i-k\right)  =k\text{)}}}\underbrace{\dbinom
{b}{a+b-\left(  i-k\right)  }}_{\substack{=\dbinom{b}{\left(  a+b\right)
+k-i}\\\text{(since }a+b-\left(  i-k\right)  =\left(  a+b\right)
+k-i\text{)}}}\nonumber\\
&  \ \ \ \ \ \ \ \ \ \ \left(  \text{here, we have substituted }i-k\text{ for
}j\text{ in the sum}\right) \nonumber\\
&  =\sum_{k=0}^{i}\dbinom{i-b}{k}\dbinom{b}{\left(  a+b\right)  +k-i}.
\label{pf.prop.binom.bin-id.f.8.pf.1}%
\end{align}
On the other hand, we have $b\in\mathbb{N}$, $\left(  i-b\right)  +b=i\geq0$
and $a\geq i-b$ (since $a+b\geq i$). Therefore, we can apply Proposition
\ref{prop.vandermonde.consequences} \textbf{(g)} to $i-b$, $b$ and $a$ instead
of $x$, $y$ and $n$. As a result, we obtain%
\begin{align*}
\dbinom{\left(  i-b\right)  +b}{a}  &  =\sum_{k=0}^{\left(  i-b\right)
+b}\dbinom{i-b}{k}\underbrace{\dbinom{b}{a+k-\left(  i-b\right)  }%
}_{\substack{=\dbinom{b}{\left(  a+b\right)  +k-i}\\\text{(since }a+k-\left(
i-b\right)  =\left(  a+b\right)  +k-i\text{)}}}\\
&  =\sum_{k=0}^{\left(  i-b\right)  +b}\dbinom{i-b}{k}\dbinom{b}{\left(
a+b\right)  +k-i}.
\end{align*}
Since $\left(  i-b\right)  +b=i$, this rewrites as
\[
\dbinom{i}{a}=\sum_{k=0}^{i}\dbinom{i-b}{k}\dbinom{b}{\left(  a+b\right)
+k-i}.
\]
Compared with (\ref{pf.prop.binom.bin-id.f.8.pf.1}), this yields%
\[
\sum_{j=0}^{i}\dbinom{i-b}{i-j}\dbinom{b}{a+b-j}=\dbinom{i}{a}.
\]
This proves (\ref{pf.prop.binom.bin-id.f.8}).]

Hence, (\ref{pf.prop.binom.bin-id.f.7}) becomes%
\begin{align*}
&  \sum_{i=0}^{b}\dbinom{a}{i}\dbinom{b}{i}\dbinom{x+i}{a+b}\\
&  =\sum_{i=0}^{a+b}\dbinom{x}{i}\dbinom{a}{a+b-i}\underbrace{\sum_{j=0}%
^{i}\dbinom{i-b}{i-j}\dbinom{b}{a+b-j}}_{\substack{=\dbinom{i}{a}\\\text{(by
(\ref{pf.prop.binom.bin-id.f.8}))}}}\\
&  =\sum_{i=0}^{a+b}\underbrace{\dbinom{x}{i}\dbinom{a}{a+b-i}\dbinom{i}{a}%
}_{=\dbinom{i}{a}\dbinom{a}{a+b-i}\dbinom{x}{i}}=\sum_{i=0}^{a+b}\dbinom{i}%
{a}\dbinom{a}{a+b-i}\dbinom{x}{i}\\
&  =\sum_{i=0}^{a-1}\underbrace{\dbinom{i}{a}}_{\substack{=0\\\text{(by
(\ref{eq.binom.0}), applied to }i\text{ and }a\\\text{instead of }m\text{ and
}n\text{ (since }i<a\text{))}}}\dbinom{a}{a+b-i}\dbinom{x}{i}+\sum_{i=a}%
^{a+b}\dbinom{i}{a}\dbinom{a}{a+b-i}\dbinom{x}{i}\\
&  \ \ \ \ \ \ \ \ \ \ \left(  \text{since }0\leq a\leq a+b\right) \\
&  =\underbrace{\sum_{i=0}^{a-1}0\dbinom{a}{a+b-i}\dbinom{x}{i}}_{=0}%
+\sum_{i=a}^{a+b}\dbinom{i}{a}\dbinom{a}{a+b-i}\dbinom{x}{i}\\
&  =\sum_{i=a}^{a+b}\dbinom{i}{a}\dbinom{a}{a+b-i}\dbinom{x}{i}=\dbinom{x}%
{a}\dbinom{x}{b}\ \ \ \ \ \ \ \ \ \ \left(  \text{by
(\ref{pf.prop.binom.bin-id.f.1})}\right)  .
\end{align*}
This proves Proposition \ref{prop.binom.bin-id} \textbf{(f)}.

\textbf{(g)} Let $a\in\mathbb{N}$, $b\in\mathbb{N}$ and $x\in\mathbb{Z}$. From
(\ref{eq.binom.upper-neg}) (applied to $m=-x-1$ and $n=a$), we obtain
$\dbinom{-x-1}{a}=\left(  -1\right)  ^{a}\dbinom{a-\left(  -x-1\right)  -1}%
{a}=\left(  -1\right)  ^{a}\dbinom{a+x}{a}$ (since $a-\left(  -x-1\right)
-1=a+x$). The same argument (applied to $b$ instead of $a$) shows that
$\dbinom{-x-1}{b}=\left(  -1\right)  ^{b}\dbinom{b+x}{b}$.

Now, Proposition \ref{prop.binom.bin-id} \textbf{(f)} (applied to $-x-1$
instead of $x$) shows that%
\begin{align}
\sum_{i=0}^{b}\dbinom{a}{i}\dbinom{b}{i}\dbinom{\left(  -x-1\right)  +i}{a+b}
&  =\underbrace{\dbinom{-x-1}{a}}_{=\left(  -1\right)  ^{a}\dbinom{a+x}{a}%
}\underbrace{\dbinom{-x-1}{b}}_{=\left(  -1\right)  ^{b}\dbinom{b+x}{b}%
}\nonumber\\
&  =\left(  -1\right)  ^{a}\dbinom{a+x}{a}\left(  -1\right)  ^{b}\dbinom
{b+x}{b}\nonumber\\
&  =\underbrace{\left(  -1\right)  ^{a}\left(  -1\right)  ^{b}}_{=\left(
-1\right)  ^{a+b}}\dbinom{a+x}{a}\dbinom{b+x}{b}\nonumber\\
&  =\left(  -1\right)  ^{a+b}\dbinom{a+x}{a}\dbinom{b+x}{b}.
\label{pf.prop.binom.bin-id.g.3}%
\end{align}
But every $i\in\left\{  0,1,\ldots,b\right\}  $ satisfies%
\begin{align*}
\dbinom{\left(  -x-1\right)  +i}{a+b}  &  =\left(  -1\right)  ^{a+b}%
\dbinom{a+b-\left(  \left(  -x-1\right)  +i\right)  -1}{a+b}\\
&  \ \ \ \ \ \ \ \ \ \ \left(  \text{by (\ref{eq.binom.upper-neg}), applied to
}m=\left(  -x-1\right)  +i\text{ and }n=a+b\right) \\
&  =\left(  -1\right)  ^{a+b}\dbinom{a+b+x-i}{a+b}\\
&  \ \ \ \ \ \ \ \ \ \ \left(  \text{since }a+b-\left(  \left(  -x-1\right)
+i\right)  -1=a+b+x-i\right)  .
\end{align*}
Hence,%
\begin{align*}
&  \sum_{i=0}^{b}\dbinom{a}{i}\dbinom{b}{i}\underbrace{\dbinom{\left(
-x-1\right)  +i}{a+b}}_{=\left(  -1\right)  ^{a+b}\dbinom{a+b+x-i}{a+b}}\\
&  =\sum_{i=0}^{b}\dbinom{a}{i}\dbinom{b}{i}\left(  -1\right)  ^{a+b}%
\dbinom{a+b+x-i}{a+b}=\left(  -1\right)  ^{a+b}\sum_{i=0}^{b}\dbinom{a}%
{i}\dbinom{b}{i}\dbinom{a+b+x-i}{a+b}.
\end{align*}
Comparing this with (\ref{pf.prop.binom.bin-id.g.3}), we obtain%
\[
\left(  -1\right)  ^{a+b}\sum_{i=0}^{b}\dbinom{a}{i}\dbinom{b}{i}%
\dbinom{a+b+x-i}{a+b}=\left(  -1\right)  ^{a+b}\dbinom{a+x}{a}\dbinom{b+x}%
{b}.
\]
We can cancel $\left(  -1\right)  ^{a+b}$ from this equality (since $\left(
-1\right)  ^{a+b}\neq0$), and thus obtain $\sum_{i=0}^{b}\dbinom{a}{i}%
\dbinom{b}{i}\dbinom{a+b+x-i}{a+b}=\dbinom{a+x}{a}\dbinom{b+x}{b}$. This
proves Proposition \ref{prop.binom.bin-id} \textbf{(g)}.
\end{proof}

Many more examples of equalities with binomial coefficients, as well as
advanced tactics for proving such equalities, can be found in \cite[Chapter
5]{GKP}.

\begin{exercise}
\label{exe.binom.rewr-prod}Let $n\in\mathbb{Q}$, $a\in\mathbb{N}$ and
$b\in\mathbb{N}$.

\textbf{(a)} Prove that every integer $j\geq a$ satisfies%
\[
\dbinom{n}{j}\dbinom{j}{a}\dbinom{n-j}{b}=\dbinom{n}{a}\dbinom{n-a}{b}%
\dbinom{n-a-b}{j-a}.
\]

\textbf{(b)} Compute the sum $\sum_{j=a}^{n}\dbinom{n}{j}\dbinom{j}{a}%
\dbinom{n-j}{b}$ for every integer $n\geq a$. (The result should contain no
summation signs.)
\end{exercise}

\subsection{The principle of inclusion and exclusion}

We shall next discuss the \textit{principle of inclusion and exclusion}, and
some of its generalizations. This is a crucial result in combinatorics, which
can help both in answering enumerative questions (i.e., questions of the form
\textquotedblleft how many objects of a given kind satisfy a certain set of
properties\textquotedblright) and in proving combinatorial identities (such
as, to give a simple example, Proposition \ref{prop.binom.bin-id}
\textbf{(c)}, but also various deeper results). We shall not dwell on the
applications of this principle; the reader can easily find them in textbooks
on enumerative combinatorics (such as \cite[\S 5.1]{Aigner07} or
\cite[\S 16]{Galvin} or \cite[Chapter 4]{Loehr-BC} or \cite[Chapter
IV]{Comtet74} or \cite[\S 15.9]{LeLeMe16} or \cite[Chapter 6]{AndFen04}). We
will, however, prove the principle and a few of its generalizations.

The principle itself (in one of its most basic forms) answers the following
simple question: Given $n$ finite sets $A_{1},A_{2},\ldots,A_{n}$, how do we
compute the size of their union $A_{1}\cup A_{2}\cup\cdots\cup A_{n}$ if we
know the sizes of all of their intersections (not just the intersection
$A_{1}\cap A_{2}\cap\cdots\cap A_{n}$, but also the intersections of some of
the sets only)? Let us first answer this question for specific small values of
$n$:

\begin{itemize}
\item For $n=1$, we have the tautological equality%
\begin{equation}
\left\vert A_{1}\right\vert =\left\vert A_{1}\right\vert . \label{eq.pie.n=1}%
\end{equation}
(Only a true mathematician would begin a study with such a statement.)

\item For $n=2$, we have the known formula%
\begin{equation}
\left\vert A_{1}\cup A_{2}\right\vert =\left\vert A_{1}\right\vert +\left\vert
A_{2}\right\vert -\left\vert A_{1}\cap A_{2}\right\vert . \label{eq.pie.n=2}%
\end{equation}
Notice that $\left\vert A_{1}\right\vert $ and $\left\vert A_{2}\right\vert $
are sizes of intersections of some of the sets $A_{1},A_{2}$: Namely, $A_{1}$
is the intersection of the single set $A_{1}$, while $A_{2}$ is the
intersection of the single set $A_{2}$.

\item For $n=3$, we have%
\begin{align}
\left\vert A_{1}\cup A_{2}\cup A_{3}\right\vert  &  =\left\vert A_{1}%
\right\vert +\left\vert A_{2}\right\vert +\left\vert A_{3}\right\vert
-\left\vert A_{1}\cap A_{2}\right\vert -\left\vert A_{1}\cap A_{3}\right\vert
-\left\vert A_{2}\cap A_{3}\right\vert \nonumber\\
&  \ \ \ \ \ \ \ \ \ \ +\left\vert A_{1}\cap A_{2}\cap A_{3}\right\vert .
\label{eq.pie.n=3}%
\end{align}
This is not as well-known as (\ref{eq.pie.n=2}), but can be easily derived by
applying (\ref{eq.pie.n=2}) twice. (In fact, first apply (\ref{eq.pie.n=2}) to
$A_{1}\cup A_{2}$ and $A_{3}$ instead of $A_{1}$ and $A_{2}$; then, apply
(\ref{eq.pie.n=2}) directly to rewrite $\left\vert A_{1}\cup A_{2}\right\vert
$.)

\item For $n=4$, we have%
\begin{align}
&  \left\vert A_{1}\cup A_{2}\cup A_{3}\cup A_{4}\right\vert \nonumber\\
&  =\left\vert A_{1}\right\vert +\left\vert A_{2}\right\vert +\left\vert
A_{3}\right\vert +\left\vert A_{4}\right\vert \nonumber\\
&  \ \ \ \ \ \ \ \ \ \ -\left\vert A_{1}\cap A_{2}\right\vert -\left\vert
A_{1}\cap A_{3}\right\vert -\left\vert A_{1}\cap A_{4}\right\vert -\left\vert
A_{2}\cap A_{3}\right\vert -\left\vert A_{2}\cap A_{4}\right\vert -\left\vert
A_{3}\cap A_{4}\right\vert \nonumber\\
&  \ \ \ \ \ \ \ \ \ \ +\left\vert A_{1}\cap A_{2}\cap A_{3}\right\vert
+\left\vert A_{1}\cap A_{2}\cap A_{4}\right\vert +\left\vert A_{1}\cap
A_{3}\cap A_{4}\right\vert +\left\vert A_{2}\cap A_{3}\cap A_{4}\right\vert
\nonumber\\
&  \ \ \ \ \ \ \ \ \ \ -\left\vert A_{1}\cap A_{2}\cap A_{3}\cap
A_{4}\right\vert . \label{eq.pie.n=4}%
\end{align}
Again, this can be derived by applying (\ref{eq.pie.n=2}) many times.
\end{itemize}

The four equalities (\ref{eq.pie.n=1}), (\ref{eq.pie.n=2}), (\ref{eq.pie.n=3})
and (\ref{eq.pie.n=4}) all follow the same pattern: On the left hand side is
the size $\left\vert A_{1}\cup A_{2}\cup\cdots\cup A_{n}\right\vert $ of the
union $A_{1}\cup A_{2}\cup\cdots\cup A_{n}$ of all the $n$ sets $A_{1}%
,A_{2},\ldots,A_{n}$. On the right hand side is an \textquotedblleft
alternating sum\textquotedblright\ (i.e., a sum, but with minus signs in front
of some of its addends), whose addends are the sizes of the intersections of
all possible choices of \textbf{some} of the $n$ sets $A_{1},A_{2}%
,\ldots,A_{n}$ (except for the choice where none of the $n$ sets are chosen;
this does not have a well-defined intersection). Notice that there are
$2^{n}-1$ such choices, so the right hand side is an \textquotedblleft
alternating sum\textquotedblright\ of $2^{n}-1$ addends. In other words, each
addend on the right hand side has the form $\left\vert A_{i_{1}}\cap A_{i_{2}%
}\cap\cdots\cap A_{i_{m}}\right\vert $ for some $m$-tuple $\left(  i_{1}%
,i_{2},\ldots,i_{m}\right)  $ of integers between $1$ and $n$ (inclusive) such
that $m\geq1$ and $i_{1}<i_{2}<\cdots<i_{m}$. The sign in front of this addend
is a plus sign if $m$ is odd, and is a minus sign if $m$ is even. Thus, we can
replace this sign by a factor of $\left(  -1\right)  ^{m-1}$.

We can try and generalize the pattern as follows:%
\begin{align}
&  \left\vert A_{1}\cup A_{2}\cup\cdots\cup A_{n}\right\vert \nonumber\\
&  =\sum_{m=1}^{n}\sum_{1\leq i_{1}<i_{2}<\cdots<i_{m}\leq n}\left(
-1\right)  ^{m-1}\left\vert A_{i_{1}}\cap A_{i_{2}}\cap\cdots\cap A_{i_{m}%
}\right\vert . \label{eq.pie.n=n-basic}%
\end{align}
Here, the summation sign \textquotedblleft$\sum_{1\leq i_{1}<i_{2}%
<\cdots<i_{m}\leq n}$\textquotedblright\ on the right hand side is an
abbreviation for $\sum_{\substack{\left(  i_{1},i_{2},\ldots,i_{m}\right)
\in\left\{  1,2,\ldots,n\right\}  ^{m};\\i_{1}<i_{2}<\cdots<i_{m}}}$.

The equality (\ref{eq.pie.n=n-basic}) is indeed correct; it is one of several
(equivalent) versions of the principle of inclusion and exclusion. For
example, it appears in \cite[\S 4.7]{Loehr-BC}. We shall, however, state it
differently, for the sake of better generalizability. First, we will index the
intersections of \textbf{some} of the $n$ sets $A_{1},A_{2},\ldots,A_{n}$ not
by $m$-tuples $\left(  i_{1},i_{2},\ldots,i_{m}\right)  \in\left\{
1,2,\ldots,n\right\}  ^{m}$ satisfying $i_{1}<i_{2}<\cdots<i_{m}$, but rather
by nonempty subsets of $\left\{  1,2,\ldots,n\right\}  $. Second, our sets
$A_{1},A_{2},\ldots,A_{n}$ will be labelled not by the numbers $1,2,\ldots,n$,
but rather by elements of a finite set $G$. This will result in a more
abstract, but also more flexible version of (\ref{eq.pie.n=n-basic}).

First, we introduce some notations:

\begin{definition}
Let $I$ be a set. For each $i\in I$, we let $A_{i}$ be a set.

\textbf{(a)} Then, $\bigcup_{i\in I}A_{i}$ denotes the union of all the sets
$A_{i}$ for $i\in I$. This union is defined by%
\[
\bigcup_{i\in I}A_{i}=\left\{  x\ \mid\ \text{there exists an }i\in I\text{
such that }x\in A_{i}\right\}  .
\]

For example, if $I=\left\{  i_{1},i_{2},\ldots,i_{k}\right\}  $ is a finite
set, then
\[
\bigcup_{i\in I}A_{i}=A_{i_{1}}\cup A_{i_{2}}\cup\cdots\cup A_{i_{k}}.
\]

Notice that $A_{j}\subseteq\bigcup_{i\in I}A_{i}$ for each $j\in I$. If $I$ is
finite, and if each of the sets $A_{i}$ is finite, then their union
$\bigcup_{i\in I}A_{i}$ is also finite.

Note that $\bigcup_{i\in\varnothing}A_{i}=\varnothing$.

\textbf{(b)} Assume that $I$ is nonempty. Then, $\bigcap_{i\in I}A_{i}$
denotes the intersection of all the sets $A_{i}$ for $i\in I$. This
intersection is defined by%
\[
\bigcap_{i\in I}A_{i}=\left\{  x\ \mid\ x\in A_{i}\text{ for all }i\in
I\right\}  .
\]

For example, if $I=\left\{  i_{1},i_{2},\ldots,i_{k}\right\}  $ is a finite
set, then
\[
\bigcap_{i\in I}A_{i}=A_{i_{1}}\cap A_{i_{2}}\cap\cdots\cap A_{i_{k}}.
\]

Notice that $\bigcap_{i\in I}A_{i}\subseteq A_{j}$ for each $j\in I$. If each
of the sets $A_{i}$ is finite, then their intersection $\bigcap_{i\in I}A_{i}$
is also finite.

\textbf{Caution:} The intersection $\bigcap_{i\in I}A_{i}$ is not defined when
$I$ is empty, because this intersection would have to contain every object in
the universe (which is impossible for a set).
\end{definition}

\begin{definition}
If $G$ is any finite set, then the sign $\sum_{I\subseteq G}$ shall mean
$\sum_{I\in\mathcal{P}\left(  G\right)  }$, where $\mathcal{P}\left(
G\right)  $ denotes the powerset of $G$. For example,%
\begin{align*}
\sum_{I\subseteq\left\{  7,8\right\}  }\prod_{i\in I}i  &  =\sum
_{I\in\mathcal{P}\left(  \left\{  7,8\right\}  \right)  }\prod_{i\in
I}i=\underbrace{\prod_{i\in\varnothing}i}_{=\left(  \text{empty product}%
\right)  =1}+\underbrace{\prod_{i\in\left\{  7\right\}  }i}_{=7}%
+\underbrace{\prod_{i\in\left\{  8\right\}  }i}_{=8}+\underbrace{\prod
_{i\in\left\{  7,8\right\}  }i}_{=7\cdot8}\\
&  \ \ \ \ \ \ \ \ \ \ \left(  \text{since the subsets of }\left\{
7,8\right\}  \text{ are }\varnothing,\left\{  7\right\}  ,\left\{  8\right\}
,\left\{  7,8\right\}  \right) \\
&  =1+7+8+7\cdot8.
\end{align*}

\end{definition}

We are now ready to state one of the forms of the principle:

\begin{theorem}
\label{thm.pie.union}Let $G$ be a finite set. For each $i\in G$, let $A_{i}$
be a finite set. Then,%
\[
\left\vert \bigcup_{i\in G}A_{i}\right\vert =\sum_{\substack{I\subseteq
G;\\I\neq\varnothing}}\left(  -1\right)  ^{\left\vert I\right\vert
-1}\left\vert \bigcap_{i\in I}A_{i}\right\vert .
\]

\end{theorem}

If $G=\left\{  1,2,\ldots,n\right\}  $ for some $n\in\mathbb{N}$, then the
formula in Theorem \ref{thm.pie.union} is a restatement of
(\ref{eq.pie.n=n-basic}) (because the nonempty subsets of $\left\{
1,2,\ldots,n\right\}  $ are in 1-to-1 correspondence with the $m$-tuples
$\left(  i_{1},i_{2},\ldots,i_{m}\right)  \in\left\{  1,2,\ldots,n\right\}
^{m}$ satisfying $i_{1}<i_{2}<\cdots<i_{m}$ and $m\in\left\{  1,2,\ldots
,n\right\}  $).

A statement equivalent to Theorem \ref{thm.pie.union} is the following:

\begin{theorem}
\label{thm.pie.nonunion}Let $S$ be a finite set. Let $G$ be a finite set. For
each $i\in G$, let $A_{i}$ be a subset of $S$. We define the intersection
$\bigcap_{i\in\varnothing}A_{i}$ (which would otherwise be undefined, since
$\varnothing$ is the empty set) to mean the set $S$. (Thus, $\bigcap_{i\in
I}A_{i}$ is defined for any subset $I$ of $G$, not just for nonempty subsets
$I$.) Then,%
\[
\left\vert S\setminus\left(  \bigcup_{i\in G}A_{i}\right)  \right\vert
=\sum_{I\subseteq G}\left(  -1\right)  ^{\left\vert I\right\vert }\left\vert
\bigcap_{i\in I}A_{i}\right\vert .
\]

\end{theorem}

Theorem \ref{thm.pie.union} and Theorem \ref{thm.pie.nonunion} are commonly
known as the \textit{principle of inclusion and exclusion}, or as the
\textit{Sylvester sieve formula}. They are not hard to prove (see, e.g.,
\cite[\S 16]{Galvin} for two proofs, and \cite[proof of (2.5)]{Sagan19} for
yet another). Rather than proving them directly, we shall however generalize
them and prove the generalization, from which they are easily obtained as
particular cases.

We generalize these theorems in two steps. First, we observe that the
$S\setminus\left(  \bigcup_{i\in G}A_{i}\right)  $ in Theorem
\ref{thm.pie.nonunion} is simply the set of all elements of $S$ that belong to
none of the subsets $A_{i}$ (for $i\in G$). In other words,%
\[
S\setminus\left(  \bigcup_{i\in G}A_{i}\right)  =\left\{  s\in S\ \mid
\ \text{the number of }i\in G\text{ satisfying }s\in A_{i}\text{ equals
}0\right\}  .
\]
We can replace the \textquotedblleft$0$\textquotedblright\ here by any number
$k$, and ask for the size of the resulting set. The answer is given by the
following result of Charles Jordan (see \cite[\S 4.8, Theorem A]{Comtet74} and
\cite{DanRot78} for fairly complicated proofs):

\begin{theorem}
\label{thm.pie.jordan}Let $S$ be a finite set. Let $G$ be a finite set. For
each $i\in G$, let $A_{i}$ be a subset of $S$. We define the intersection
$\bigcap_{i\in\varnothing}A_{i}$ (which would otherwise be undefined, since
$\varnothing$ is the empty set) to mean the set $S$. (Thus, $\bigcap_{i\in
I}A_{i}$ is defined for any subset $I$ of $G$, not just for nonempty subsets
$I$.)

Let $k\in\mathbb{N}$. Let
\[
S_{k}=\left\{  s\in S\ \mid\ \text{the number of }i\in G\text{ satisfying
}s\in A_{i}\text{ equals }k\right\}  .
\]
(In other words, $S_{k}$ is the set of all elements of $S$ that belong to
exactly $k$ of the subsets $A_{i}$.) Then,%
\[
\left\vert S_{k}\right\vert =\sum_{I\subseteq G}\left(  -1\right)
^{\left\vert I\right\vert -k}\dbinom{\left\vert I\right\vert }{k}\left\vert
\bigcap_{i\in I}A_{i}\right\vert .
\]

\end{theorem}

A different generalization of Theorem \ref{thm.pie.nonunion} (closely related
to the \textit{Bonferroni inequalities}, for which see \cite[\S 17, problem
(2)]{Galvin}) explores what happens when the sum on the right hand side of the
formula is restricted to only those subsets $I$ of $G$ whose size doesn't
surpass a given integer $m$:

\begin{theorem}
\label{thm.pie.bonfer}Let $S$ be a finite set. Let $G$ be a finite set. For
each $i\in G$, let $A_{i}$ be a subset of $S$. We define the intersection
$\bigcap_{i\in\varnothing}A_{i}$ (which would otherwise be undefined, since
$\varnothing$ is the empty set) to mean the set $S$. (Thus, $\bigcap_{i\in
I}A_{i}$ is defined for any subset $I$ of $G$, not just for nonempty subsets
$I$.)

Let $m\in\mathbb{N}$. For each $s\in S$, let $c\left(  s\right)  $ denote the
number of $i\in G$ satisfying $s\in A_{i}$. Then,%
\[
\left(  -1\right)  ^{m}\sum_{s\in S}\dbinom{c\left(  s\right)  -1}{m}%
=\sum_{\substack{I\subseteq G;\\\left\vert I\right\vert \leq m}}\left(
-1\right)  ^{\left\vert I\right\vert }\left\vert \bigcap_{i\in I}%
A_{i}\right\vert .
\]

\end{theorem}

Finally, Theorem \ref{thm.pie.jordan} and Theorem \ref{thm.pie.bonfer} can be
merged into one common general principle:

\begin{theorem}
\label{thm.pie.merged}Let $S$ be a finite set. Let $G$ be a finite set. For
each $i\in G$, let $A_{i}$ be a subset of $S$. We define the intersection
$\bigcap_{i\in\varnothing}A_{i}$ (which would otherwise be undefined, since
$\varnothing$ is the empty set) to mean the set $S$. (Thus, $\bigcap_{i\in
I}A_{i}$ is defined for any subset $I$ of $G$, not just for nonempty subsets
$I$.)

Let $k\in\mathbb{N}$ and $m\in\mathbb{N}$ be such that $m\geq k$. For each
$s\in S$, let $c\left(  s\right)  $ denote the number of $i\in G$ satisfying
$s\in A_{i}$. Then,%
\[
\left(  -1\right)  ^{m}\sum_{s\in S}\dbinom{c\left(  s\right)  }{k}%
\dbinom{c\left(  s\right)  -k-1}{m-k}=\sum_{\substack{I\subseteq
G;\\\left\vert I\right\vert \leq m}}\left(  -1\right)  ^{\left\vert
I\right\vert }\dbinom{\left\vert I\right\vert }{k}\left\vert \bigcap_{i\in
I}A_{i}\right\vert .
\]

\end{theorem}

As we said, we shall first prove Theorem \ref{thm.pie.merged}, and then derive
all the preceding theorems in this section from it. The proof of Theorem
\ref{thm.pie.merged} will rely on several ingredients, the first of which is
the following simple identity:

\begin{lemma}
\label{lem.pie.binid}Let $n\in\mathbb{N}$ and $k\in\mathbb{N}$. Let
$m\in\left\{  k,k+1,k+2,\ldots\right\}  $. Then,%
\[
\sum_{r=0}^{m}\left(  -1\right)  ^{r}\dbinom{n}{r}\dbinom{r}{k}=\left(
-1\right)  ^{m}\dbinom{n}{k}\dbinom{n-k-1}{m-k}.
\]

\end{lemma}

\begin{exercise}
\label{exe.pie.binid}Prove Lemma \ref{lem.pie.binid}.
\end{exercise}

Next, we introduce a simple yet immensely helpful notation that will
facilitate our proof:

\begin{definition}
\label{def.iverson}If $\mathcal{A}$ is any logical statement, then we define
an integer $\left[  \mathcal{A}\right]  \in\left\{  0,1\right\}  $ by%
\[
\left[  \mathcal{A}\right]  =%
\begin{cases}
1, & \text{if }\mathcal{A}\text{ is true};\\
0, & \text{if }\mathcal{A}\text{ is false}%
\end{cases}
.
\]
For example, $\left[  1+1=2\right]  =1$ (since $1+1=2$ is true), whereas
$\left[  1+1=1\right]  =0$ (since $1+1=1$ is false).

If $\mathcal{A}$ is any logical statement, then the integer $\left[
\mathcal{A}\right]  $ is known as the \textit{truth value} of $\mathcal{A}$.
The notation $\left[  \mathcal{A}\right]  $ is known as
\href{https://en.wikipedia.org/wiki/Iverson_bracket}{the \textit{Iverson
bracket notation}}.
\end{definition}

Clearly, if $\mathcal{A}$ and $\mathcal{B}$ are two equivalent logical
statements, then $\left[  \mathcal{A}\right]  =\left[  \mathcal{B}\right]  $.
This and a few other useful properties of the Iverson bracket notation are
collected in the following exercise:

\begin{exercise}
\label{exe.iverson-prop}Prove the following rules for truth values:

\textbf{(a)} If $\mathcal{A}$ and $\mathcal{B}$ are two equivalent logical
statements, then $\left[  \mathcal{A}\right]  =\left[  \mathcal{B}\right]  $.

\textbf{(b)} If $\mathcal{A}$ is any logical statement, then $\left[
\text{not }\mathcal{A}\right]  =1-\left[  \mathcal{A}\right]  $.

\textbf{(c)} If $\mathcal{A}$ and $\mathcal{B}$ are two logical statements,
then $\left[  \mathcal{A}\wedge\mathcal{B}\right]  =\left[  \mathcal{A}%
\right]  \left[  \mathcal{B}\right]  $.

\textbf{(d)} If $\mathcal{A}$ and $\mathcal{B}$ are two logical statements,
then $\left[  \mathcal{A}\vee\mathcal{B}\right]  =\left[  \mathcal{A}\right]
+\left[  \mathcal{B}\right]  -\left[  \mathcal{A}\right]  \left[
\mathcal{B}\right]  $.

\textbf{(e)} If $\mathcal{A}$, $\mathcal{B}$ and $\mathcal{C}$ are three
logical statements, then%
\[
\left[  \mathcal{A}\vee\mathcal{B}\vee\mathcal{C}\right]  =\left[
\mathcal{A}\right]  +\left[  \mathcal{B}\right]  +\left[  \mathcal{C}\right]
-\left[  \mathcal{A}\right]  \left[  \mathcal{B}\right]  -\left[
\mathcal{A}\right]  \left[  \mathcal{C}\right]  -\left[  \mathcal{B}\right]
\left[  \mathcal{C}\right]  +\left[  \mathcal{A}\right]  \left[
\mathcal{B}\right]  \left[  \mathcal{C}\right]  .
\]

\end{exercise}

The Iverson bracket helps us rewrite the cardinality of a set as a sum:

\begin{lemma}
\label{lem.iverson.card}Let $S$ be a finite set. Let $T$ be a subset of $S$.
Then,%
\[
\left\vert T\right\vert =\sum_{s\in S}\left[  s\in T\right]  .
\]

\end{lemma}

\begin{proof}
[Proof of Lemma \ref{lem.iverson.card}.]We have%
\begin{align*}
\sum_{s\in S}\left[  s\in T\right]   &  =\underbrace{\sum_{\substack{s\in
S;\\s\in T}}}_{\substack{=\sum_{s\in T}\\\text{(since }T\text{ is
a}\\\text{subset of }S\text{)}}}\underbrace{\left[  s\in T\right]
}_{\substack{=1\\\text{(since }s\in T\text{ is true)}}}+\sum_{\substack{s\in
S;\\s\notin T}}\underbrace{\left[  s\in T\right]  }%
_{\substack{=0\\\text{(since }s\in T\text{ is false}\\\text{(since }s\notin
T\text{))}}}\\
&  \ \ \ \ \ \ \ \ \ \ \ \ \ \ \ \ \ \ \ \ \left(
\begin{array}
[c]{c}%
\text{since each }s\in S\text{ satisfies}\\
\text{either }s\in T\text{ or }s\notin T\text{ (but not both)}%
\end{array}
\right) \\
&  =\underbrace{\sum_{s\in T}1}_{=\left\vert T\right\vert \cdot1}%
+\underbrace{\sum_{\substack{s\in S;\\s\notin T}}0}_{=0}=\left\vert
T\right\vert \cdot1+0=\left\vert T\right\vert .
\end{align*}
This proves Lemma \ref{lem.iverson.card}.
\end{proof}

We can now state the main precursor to Theorem \ref{thm.pie.merged}:

\begin{lemma}
\label{lem.pie.merged.iverson}Let $S$ be a finite set. Let $G$ be a finite
set. For each $i\in G$, let $A_{i}$ be a subset of $S$. We define the
intersection $\bigcap_{i\in\varnothing}A_{i}$ (which would otherwise be
undefined, since $\varnothing$ is the empty set) to mean the set $S$. (Thus,
$\bigcap_{i\in I}A_{i}$ is defined for any subset $I$ of $G$, not just for
nonempty subsets $I$.)

Let $k\in\mathbb{N}$ and $m\in\mathbb{N}$ be such that $m\geq k$.

Let $s\in S$. Let $c\left(  s\right)  $ denote the number of $i\in G$
satisfying $s\in A_{i}$. Then,%
\[
\sum_{\substack{I\subseteq G;\\\left\vert I\right\vert \leq m}}\left(
-1\right)  ^{\left\vert I\right\vert }\dbinom{\left\vert I\right\vert }%
{k}\left[  s\in\bigcap_{i\in I}A_{i}\right]  =\left(  -1\right)  ^{m}%
\dbinom{c\left(  s\right)  }{k}\dbinom{c\left(  s\right)  -k-1}{m-k}.
\]

\end{lemma}

\begin{vershort}
\begin{proof}
[Proof of Lemma \ref{lem.pie.merged.iverson}.]From $m\geq k$ and
$m\in\mathbb{N}$, we obtain $m\in\left\{  k,k+1,k+2,\ldots\right\}  $.

Define a subset $C$ of $G$ by%
\[
C=\left\{  i\in G\ \mid\ s\in A_{i}\right\}  .
\]
Thus,%
\begin{align*}
\left\vert C\right\vert  &  =\left\vert \left\{  i\in G\ \mid\ s\in
A_{i}\right\}  \right\vert \\
&  =\left(  \text{the number of }i\in G\text{ satisfying }s\in A_{i}\right)
=c\left(  s\right)
\end{align*}
(since $c\left(  s\right)  $ was defined as the number of $i\in G$ satisfying
$s\in A_{i}$). In other words, $C$ is a $c\left(  s\right)  $-element set.
Hence, for each $r\in\mathbb{N}$, Proposition \ref{prop.binom.subsets}
(applied to $c\left(  s\right)  $, $r$ and $C$ instead of $m$, $n$ and $S$)
shows that
\begin{equation}
\dbinom{c\left(  s\right)  }{r}\text{ is the number of all }r\text{-element
subsets of }C. \label{pf.lem.pie.merged.iverson.short.number}%
\end{equation}

Let $I$ be a subset of $G$. We have the following equivalence:%
\begin{equation}
\left(  s\in\bigcap_{i\in I}A_{i}\right)  \ \Longleftrightarrow\ \left(  s\in
A_{i}\text{ for all }i\in I\right)
\label{pf.lem.pie.merged.iverson.short.equiv1}%
\end{equation}
\footnote{\textit{Proof of (\ref{pf.lem.pie.merged.iverson.short.equiv1}):} If
$I$ is nonempty, then the equivalence
(\ref{pf.lem.pie.merged.iverson.short.equiv1}) follows immediately from the
equality
\[
\bigcap_{i\in I}A_{i}=\left\{  x\ \mid\ x\in A_{i}\text{ for all }i\in
I\right\}
\]
(which is the definition of the intersection $\bigcap_{i\in I}A_{i}$). Thus,
for the rest of this proof, we WLOG assume that $I$ is not nonempty.
\par
Hence, the set $I$ is empty. In other words, $I=\varnothing$. Hence,
$\bigcap_{i\in I}A_{i}=\bigcap_{i\in\varnothing}A_{i}=S$. Thus, $s\in
S=\bigcap_{i\in I}A_{i}$. Hence, the statement $\left(  s\in\bigcap_{i\in
I}A_{i}\right)  $ is true.
\par
Also, there exist no $i\in I$ (since the set $I$ is empty). Hence, the
statement $\left(  s\in A_{i}\text{ for all }i\in I\right)  $ is vacuously
true.
\par
Thus, the statements $\left(  s\in\bigcap_{i\in I}A_{i}\right)  $ and $\left(
s\in A_{i}\text{ for all }i\in I\right)  $ are both true, and therefore
equivalent. This proves the equivalence
(\ref{pf.lem.pie.merged.iverson.short.equiv1}).}.

But if $i\in I$, then we have the equivalence%
\begin{equation}
\left(  i\in C\right)  \ \Longleftrightarrow\ \left(  s\in A_{i}\right)
\label{pf.lem.pie.merged.iverson.short.equiv2}%
\end{equation}
(since $C=\left\{  i\in G\ \mid\ s\in A_{i}\right\}  $).

Hence, the equivalence (\ref{pf.lem.pie.merged.iverson.short.equiv1}) becomes%
\begin{align*}
\left(  s\in\bigcap_{i\in I}A_{i}\right)  \  &  \Longleftrightarrow\ \left(
\underbrace{s\in A_{i}}_{\substack{\Longleftrightarrow\ \left(  i\in C\right)
\\\text{(by (\ref{pf.lem.pie.merged.iverson.short.equiv2}))}}}\text{ for all
}i\in I\right) \\
&  \Longleftrightarrow\ \left(  i\in C\text{ for all }i\in I\right)
\ \Longleftrightarrow\ \left(  I\subseteq C\right)  .
\end{align*}
In other words, the two statements $\left(  s\in\bigcap_{i\in I}A_{i}\right)
$ and $\left(  I\subseteq C\right)  $ are equivalent. Hence, Exercise
\ref{exe.iverson-prop} \textbf{(a)} (applied to $\mathcal{A}=\left(
s\in\bigcap_{i\in I}A_{i}\right)  $ and $\mathcal{B}=\left(  I\subseteq
C\right)  $) shows that%
\begin{equation}
\left[  s\in\bigcap_{i\in I}A_{i}\right]  =\left[  I\subseteq C\right]  .
\label{pf.lem.pie.merged.iverson.short.iverson-eq}%
\end{equation}

Now, forget that we fixed $I$. We thus have proven the equality
(\ref{pf.lem.pie.merged.iverson.short.iverson-eq}) for every subset $I$ of $G$.

Now,%
\begin{align*}
&  \sum_{\substack{I\subseteq G;\\\left\vert I\right\vert \leq m}}\left(
-1\right)  ^{\left\vert I\right\vert }\dbinom{\left\vert I\right\vert }%
{k}\underbrace{\left[  s\in\bigcap_{i\in I}A_{i}\right]  }_{\substack{=\left[
I\subseteq C\right]  \\\text{(by
(\ref{pf.lem.pie.merged.iverson.short.iverson-eq}))}}}\\
&  =\sum_{\substack{I\subseteq G;\\\left\vert I\right\vert \leq m}}\left(
-1\right)  ^{\left\vert I\right\vert }\dbinom{\left\vert I\right\vert }%
{k}\left[  I\subseteq C\right] \\
&  =\underbrace{\sum_{\substack{I\subseteq G;\\\left\vert I\right\vert \leq
m;\\I\subseteq C}}}_{\substack{=\sum_{\substack{I\subseteq G;\\I\subseteq
C;\\\left\vert I\right\vert \leq m}}=\sum_{\substack{I\subseteq C;\\\left\vert
I\right\vert \leq m}}\\\text{(since }C\text{ is a subset of }G\text{)}%
}}\left(  -1\right)  ^{\left\vert I\right\vert }\dbinom{\left\vert
I\right\vert }{k}\underbrace{\left[  I\subseteq C\right]  }%
_{\substack{=1\\\text{(since }I\subseteq C\text{)}}}+\sum
_{\substack{I\subseteq G;\\\left\vert I\right\vert \leq m;\\\text{not
}I\subseteq C}}\left(  -1\right)  ^{\left\vert I\right\vert }\dbinom
{\left\vert I\right\vert }{k}\underbrace{\left[  I\subseteq C\right]
}_{\substack{=0\\\text{(since we don't have }I\subseteq C\text{)}}}\\
&  \ \ \ \ \ \ \ \ \ \ \ \ \ \ \ \ \ \ \ \ \left(
\begin{array}
[c]{c}%
\text{since each subset }I\text{ of }G\text{ satisfies either }I\subseteq C\\
\text{or }\left(  \text{not }I\subseteq C\right)  \text{ (but not both)}%
\end{array}
\right) \\
&  =\sum_{\substack{I\subseteq C;\\\left\vert I\right\vert \leq m}}\left(
-1\right)  ^{\left\vert I\right\vert }\dbinom{\left\vert I\right\vert }%
{k}+\underbrace{\sum_{\substack{I\subseteq G;\\\left\vert I\right\vert \leq
m;\\\text{not }I\subseteq C}}\left(  -1\right)  ^{\left\vert I\right\vert
}\dbinom{\left\vert I\right\vert }{k}0}_{=0}=\underbrace{\sum
_{\substack{I\subseteq C;\\\left\vert I\right\vert \leq m}}}_{=\sum_{r=0}%
^{m}\sum_{\substack{I\subseteq C;\\\left\vert I\right\vert =r}}}\left(
-1\right)  ^{\left\vert I\right\vert }\dbinom{\left\vert I\right\vert }{k}\\
&  =\sum_{r=0}^{m}\sum_{\substack{I\subseteq C;\\\left\vert I\right\vert
=r}}\underbrace{\left(  -1\right)  ^{\left\vert I\right\vert }}%
_{\substack{=\left(  -1\right)  ^{r}\\\text{(since }\left\vert I\right\vert
=r\text{)}}}\underbrace{\dbinom{\left\vert I\right\vert }{k}}%
_{\substack{=\dbinom{r}{k}\\\text{(since }\left\vert I\right\vert =r\text{)}}}
\end{align*}%
\begin{align*}
&  =\sum_{r=0}^{m}\underbrace{\sum_{\substack{I\subseteq C;\\\left\vert
I\right\vert =r}}\left(  -1\right)  ^{r}\dbinom{r}{k}}_{=\left(  \text{the
number of all subsets }I\text{ of }C\text{ satisfying }\left\vert I\right\vert
=r\right)  \left(  -1\right)  ^{r}\dbinom{r}{k}}\\
&  =\sum_{r=0}^{m}\underbrace{\left(  \text{the number of all subsets }I\text{
of }C\text{ satisfying }\left\vert I\right\vert =r\right)  }%
_{\substack{=\left(  \text{the number of all }r\text{-element subsets of
}C\right)  =\dbinom{c\left(  s\right)  }{r}\\\text{(by
(\ref{pf.lem.pie.merged.iverson.short.number}))}}}\left(  -1\right)
^{r}\dbinom{r}{k}\\
&  =\sum_{r=0}^{m}\dbinom{c\left(  s\right)  }{r}\left(  -1\right)
^{r}\dbinom{r}{k}=\sum_{r=0}^{m}\left(  -1\right)  ^{r}\dbinom{c\left(
s\right)  }{r}\dbinom{r}{k}=\left(  -1\right)  ^{m}\dbinom{c\left(  s\right)
}{k}\dbinom{c\left(  s\right)  -k-1}{m-k}\\
&  \ \ \ \ \ \ \ \ \ \ \left(  \text{by Lemma \ref{lem.pie.binid} (applied to
}n=c\left(  s\right)  \text{)}\right)  .
\end{align*}
This proves Lemma \ref{lem.pie.merged.iverson}.
\end{proof}
\end{vershort}

\begin{verlong}
\begin{proof}
[Proof of Lemma \ref{lem.pie.merged.iverson}.]From $m\geq k$ and
$m\in\mathbb{N}$, we obtain $m\in\left\{  k,k+1,k+2,\ldots\right\}  $.

Define a subset $C$ of $G$ by%
\[
C=\left\{  i\in G\ \mid\ s\in A_{i}\right\}  .
\]
Thus,%
\begin{align*}
\left\vert C\right\vert  &  =\left\vert \left\{  i\in G\ \mid\ s\in
A_{i}\right\}  \right\vert \\
&  =\left(  \text{the number of }i\in G\text{ satisfying }s\in A_{i}\right)
=c\left(  s\right)
\end{align*}
(since the number of $i\in G$ satisfying $s\in A_{i}$ is $c\left(  s\right)  $
(by the definition of $c\left(  s\right)  $)). In other words, $C$ is a
$c\left(  s\right)  $-element set.

Let $I$ be a subset of $G$. We have the following equivalence:%
\begin{equation}
\left(  s\in\bigcap_{i\in I}A_{i}\right)  \ \Longleftrightarrow\ \left(  s\in
A_{i}\text{ for all }i\in I\right)  \label{pf.lem.pie.merged.iverson.equiv1}%
\end{equation}
\footnote{\textit{Proof of (\ref{pf.lem.pie.merged.iverson.equiv1}):} If $I$
is nonempty, then the equivalence (\ref{pf.lem.pie.merged.iverson.equiv1})
follows immediately from the equality
\[
\bigcap_{i\in I}A_{i}=\left\{  x\ \mid\ x\in A_{i}\text{ for all }i\in
I\right\}
\]
(which is the definition of the intersection $\bigcap_{i\in I}A_{i}$). Thus,
for the rest of this proof, we can WLOG assume that $I$ is not nonempty.
Assume this.
\par
The set $I$ is empty (since $I$ is not nonempty). In other words,
$I=\varnothing$. Hence, $\bigcap_{i\in I}A_{i}=\bigcap_{i\in\varnothing}%
A_{i}=S$. Thus, $s\in S=\bigcap_{i\in I}A_{i}$. Hence, the statement $\left(
s\in\bigcap_{i\in I}A_{i}\right)  $ is true.
\par
Also, there exist no $i\in I$ (since the set $I$ is empty). Hence, the
statement $\left(  s\in A_{i}\text{ for all }i\in I\right)  $ is vacuously
true.
\par
Thus, the statements $\left(  s\in\bigcap_{i\in I}A_{i}\right)  $ and $\left(
s\in A_{i}\text{ for all }i\in I\right)  $ are both true. Hence, these
statements are equivalent. This proves the equivalence
(\ref{pf.lem.pie.merged.iverson.equiv1}).}.

But if $i\in I$, then we have the equivalence%
\begin{equation}
\left(  i\in C\right)  \ \Longleftrightarrow\ \left(  s\in A_{i}\right)
\label{pf.lem.pie.merged.iverson.equiv2}%
\end{equation}
(since $C=\left\{  i\in G\ \mid\ s\in A_{i}\right\}  $).

Hence, the equivalence (\ref{pf.lem.pie.merged.iverson.equiv1}) becomes%
\begin{align*}
\left(  s\in\bigcap_{i\in I}A_{i}\right)  \  &  \Longleftrightarrow\ \left(
\underbrace{s\in A_{i}}_{\substack{\Longleftrightarrow\ \left(  i\in C\right)
\\\text{(by (\ref{pf.lem.pie.merged.iverson.equiv2}))}}}\text{ for all }i\in
I\right) \\
&  \Longleftrightarrow\ \left(  i\in C\text{ for all }i\in I\right)
\ \Longleftrightarrow\ \left(  I\subseteq C\right)  .
\end{align*}
In other words, the two statements $\left(  s\in\bigcap_{i\in I}A_{i}\right)
$ and $\left(  I\subseteq C\right)  $ are equivalent. Hence, Exercise
\ref{exe.iverson-prop} \textbf{(a)} (applied to $\mathcal{A}=\left(
s\in\bigcap_{i\in I}A_{i}\right)  $ and $\mathcal{B}=\left(  I\subseteq
C\right)  $) shows that%
\begin{equation}
\left[  s\in\bigcap_{i\in I}A_{i}\right]  =\left[  I\subseteq C\right]  .
\label{pf.lem.pie.merged.iverson.iverson-eq}%
\end{equation}

Now, forget that we fixed $I$. We thus have proven the equality
(\ref{pf.lem.pie.merged.iverson.iverson-eq}) for every subset $I$ of $G$.

Also, for any $r\in\mathbb{N}$, we have%
\begin{equation}
\left(  \text{the number of all subsets }I\text{ of }C\text{ satisfying
}\left\vert I\right\vert =r\right)  =\dbinom{c\left(  s\right)  }{r}
\label{pf.lem.pie.merged.iverson.number}%
\end{equation}
\footnote{\textit{Proof of (\ref{pf.lem.pie.merged.iverson.number}):} Let
$r\in\mathbb{N}$. The set $C$ is a $c\left(  s\right)  $-element set. Thus,
Proposition \ref{prop.binom.subsets} (applied to $c\left(  s\right)  $, $r$
and $C$ instead of $m$, $n$ and $S$) yields that%
\[
\dbinom{c\left(  s\right)  }{r}\text{ is the number of all }r\text{-element
subsets of }C\text{.}%
\]
In other words,%
\begin{align*}
\dbinom{c\left(  s\right)  }{r}  &  =\left(  \text{the number of all
}r\text{-element subsets of }C\right) \\
&  =\left(  \text{the number of all subsets }I\text{ of }C\text{ such that
}\underbrace{I\text{ is an }r\text{-element set}}_{\Longleftrightarrow
\ \left(  \left\vert I\right\vert =r\right)  }\right) \\
&  =\left(  \text{the number of all subsets }I\text{ of }C\text{ such that
}\left\vert I\right\vert =r\right)  .
\end{align*}
This proves (\ref{pf.lem.pie.merged.iverson.number}).}.

Now,%
\begin{align*}
&  \sum_{\substack{I\subseteq G;\\\left\vert I\right\vert \leq m}}\left(
-1\right)  ^{\left\vert I\right\vert }\dbinom{\left\vert I\right\vert }%
{k}\underbrace{\left[  s\in\bigcap_{i\in I}A_{i}\right]  }_{\substack{=\left[
I\subseteq C\right]  \\\text{(by (\ref{pf.lem.pie.merged.iverson.iverson-eq}%
))}}}\\
&  =\sum_{\substack{I\subseteq G;\\\left\vert I\right\vert \leq m}}\left(
-1\right)  ^{\left\vert I\right\vert }\dbinom{\left\vert I\right\vert }%
{k}\left[  I\subseteq C\right] \\
&  =\underbrace{\sum_{\substack{I\subseteq G;\\\left\vert I\right\vert \leq
m;\\I\subseteq C}}}_{\substack{=\sum_{\substack{I\subseteq G;\\I\subseteq
C;\\\left\vert I\right\vert \leq m}}=\sum_{\substack{I\subseteq C;\\\left\vert
I\right\vert \leq m}}\\\text{(since }C\text{ is a subset of }G\text{)}%
}}\left(  -1\right)  ^{\left\vert I\right\vert }\dbinom{\left\vert
I\right\vert }{k}\underbrace{\left[  I\subseteq C\right]  }%
_{\substack{=1\\\text{(since }I\subseteq C\text{)}}}+\sum
_{\substack{I\subseteq G;\\\left\vert I\right\vert \leq m;\\\text{not
}I\subseteq C}}\left(  -1\right)  ^{\left\vert I\right\vert }\dbinom
{\left\vert I\right\vert }{k}\underbrace{\left[  I\subseteq C\right]
}_{\substack{=0\\\text{(since we don't have }I\subseteq C\text{)}}}\\
&  \ \ \ \ \ \ \ \ \ \ \ \ \ \ \ \ \ \ \ \ \left(
\begin{array}
[c]{c}%
\text{since each subset }I\text{ of }G\text{ satisfies either }I\subseteq C\\
\text{or }\left(  \text{not }I\subseteq C\right)  \text{ (but not both)}%
\end{array}
\right) \\
&  =\sum_{\substack{I\subseteq C;\\\left\vert I\right\vert \leq m}}\left(
-1\right)  ^{\left\vert I\right\vert }\dbinom{\left\vert I\right\vert }%
{k}+\underbrace{\sum_{\substack{I\subseteq G;\\\left\vert I\right\vert \leq
m;\\\text{not }I\subseteq C}}\left(  -1\right)  ^{\left\vert I\right\vert
}\dbinom{\left\vert I\right\vert }{k}0}_{=0}=\underbrace{\sum
_{\substack{I\subseteq C;\\\left\vert I\right\vert \leq m}}}_{\substack{=\sum
_{\substack{I\subseteq C;\\\left\vert I\right\vert \in\left\{  0,1,\ldots
,m\right\}  }}\\\text{(because for any subset }I\text{ of }C\text{,}%
\\\text{the condition }\left(  \left\vert I\right\vert \leq m\right)  \text{
is equivalent}\\\text{to }\left(  \left\vert I\right\vert \in\left\{
0,1,\ldots,m\right\}  \right)  \text{ (since }\left\vert I\right\vert
\in\mathbb{N}\text{))}}}\left(  -1\right)  ^{\left\vert I\right\vert }%
\dbinom{\left\vert I\right\vert }{k}\\
&  =\underbrace{\sum_{\substack{I\subseteq C;\\\left\vert I\right\vert
\in\left\{  0,1,\ldots,m\right\}  }}}_{=\sum_{r\in\left\{  0,1,\ldots
,m\right\}  }\sum_{\substack{I\subseteq C;\\\left\vert I\right\vert =r}%
}}\left(  -1\right)  ^{\left\vert I\right\vert }\dbinom{\left\vert
I\right\vert }{k}=\sum_{r\in\left\{  0,1,\ldots,m\right\}  }\sum
_{\substack{I\subseteq C;\\\left\vert I\right\vert =r}}\underbrace{\left(
-1\right)  ^{\left\vert I\right\vert }}_{\substack{=\left(  -1\right)
^{r}\\\text{(since }\left\vert I\right\vert =r\text{)}}}\underbrace{\dbinom
{\left\vert I\right\vert }{k}}_{\substack{=\dbinom{r}{k}\\\text{(since
}\left\vert I\right\vert =r\text{)}}}
\end{align*}%
\begin{align*}
&  =\underbrace{\sum_{r\in\left\{  0,1,\ldots,m\right\}  }}_{=\sum_{r=0}^{m}%
}\underbrace{\sum_{\substack{I\subseteq C;\\\left\vert I\right\vert
=r}}\left(  -1\right)  ^{r}\dbinom{r}{k}}_{=\left(  \text{the number of all
subsets }I\text{ of }C\text{ satisfying }\left\vert I\right\vert =r\right)
\left(  -1\right)  ^{r}\dbinom{r}{k}}\\
&  =\sum_{r=0}^{m}\underbrace{\left(  \text{the number of all subsets }I\text{
of }C\text{ satisfying }\left\vert I\right\vert =r\right)  }%
_{\substack{=\dbinom{c\left(  s\right)  }{r}\\\text{(by
(\ref{pf.lem.pie.merged.iverson.number}))}}}\left(  -1\right)  ^{r}\dbinom
{r}{k}\\
&  =\sum_{r=0}^{m}\dbinom{c\left(  s\right)  }{r}\left(  -1\right)
^{r}\dbinom{r}{k}=\sum_{r=0}^{m}\left(  -1\right)  ^{r}\dbinom{c\left(
s\right)  }{r}\dbinom{r}{k}\\
&  =\left(  -1\right)  ^{m}\dbinom{c\left(  s\right)  }{k}\dbinom{c\left(
s\right)  -k-1}{m-k}\\
&  \ \ \ \ \ \ \ \ \ \ \left(  \text{by Lemma \ref{lem.pie.binid} (applied to
}n=c\left(  s\right)  \text{)}\right)  .
\end{align*}
This proves Lemma \ref{lem.pie.merged.iverson}.
\end{proof}
\end{verlong}

We now easily obtain Theorem \ref{thm.pie.merged}:

\begin{proof}
[Proof of Theorem \ref{thm.pie.merged}.]For each subset $I$ of $G$, the
intersection $\bigcap_{i\in I}A_{i}$ is a subset of $S$%
\ \ \ \ \footnote{\textit{Proof.} Let $I$ be a subset of $G$. We must show
that the intersection $\bigcap_{i\in I}A_{i}$ is a subset of $S$.
\par
If $I$ is nonempty, then this is clear (because each of the sets $A_{i}$ is a
subset of $S$). Hence, for the rest of this proof, we can WLOG assume that $I$
is not nonempty. Assume this.
\par
The set $I$ is empty (since $I$ is not nonempty). Hence, $I=\varnothing$.
Thus, $\bigcap_{i\in I}A_{i}=\bigcap_{i\in\varnothing}A_{i}=S$ (since we
defined $\bigcap_{i\in\varnothing}A_{i}$ to be $S$). Hence, $\bigcap
_{i\in\varnothing}A_{i}$ is a subset of $S$. Qed.}. Hence, for each subset $I$
of $G$, we obtain%
\begin{equation}
\left\vert \bigcap_{i\in I}A_{i}\right\vert =\sum_{s\in S}\left[  s\in
\bigcap_{i\in I}A_{i}\right]  \label{pf.thm.pie.merged.intersection}%
\end{equation}
(by Lemma \ref{lem.iverson.card} (applied to $T=\bigcap_{i\in I}A_{i}$)).
Hence,%
\begin{align*}
&  \sum_{\substack{I\subseteq G;\\\left\vert I\right\vert \leq m}}\left(
-1\right)  ^{\left\vert I\right\vert }\dbinom{\left\vert I\right\vert }%
{k}\underbrace{\left\vert \bigcap_{i\in I}A_{i}\right\vert }_{\substack{=\sum
_{s\in S}\left[  s\in\bigcap_{i\in I}A_{i}\right]  \\\text{(by
(\ref{pf.thm.pie.merged.intersection}))}}}\\
&  =\sum_{\substack{I\subseteq G;\\\left\vert I\right\vert \leq m}}\left(
-1\right)  ^{\left\vert I\right\vert }\dbinom{\left\vert I\right\vert }{k}%
\sum_{s\in S}\left[  s\in\bigcap_{i\in I}A_{i}\right]  =\underbrace{\sum
_{\substack{I\subseteq G;\\\left\vert I\right\vert \leq m}}\sum_{s\in S}%
}_{=\sum_{s\in S}\sum_{\substack{I\subseteq G;\\\left\vert I\right\vert \leq
m}}}\left(  -1\right)  ^{\left\vert I\right\vert }\dbinom{\left\vert
I\right\vert }{k}\left[  s\in\bigcap_{i\in I}A_{i}\right] \\
&  =\sum_{s\in S}\underbrace{\sum_{\substack{I\subseteq G;\\\left\vert
I\right\vert \leq m}}\left(  -1\right)  ^{\left\vert I\right\vert }%
\dbinom{\left\vert I\right\vert }{k}\left[  s\in\bigcap_{i\in I}A_{i}\right]
}_{\substack{=\left(  -1\right)  ^{m}\dbinom{c\left(  s\right)  }{k}%
\dbinom{c\left(  s\right)  -k-1}{m-k}\\\text{(by Lemma
\ref{lem.pie.merged.iverson})}}}\\
&  =\sum_{s\in S}\left(  -1\right)  ^{m}\dbinom{c\left(  s\right)  }{k}%
\dbinom{c\left(  s\right)  -k-1}{m-k}=\left(  -1\right)  ^{m}\sum_{s\in
S}\dbinom{c\left(  s\right)  }{k}\dbinom{c\left(  s\right)  -k-1}{m-k}.
\end{align*}
This proves Theorem \ref{thm.pie.merged}.
\end{proof}

Having proven Theorem \ref{thm.pie.merged}, we can now easily derive the other
(less general) versions of the inclusion-exclusion principle:

\begin{exercise}
\label{exe.pie.specialize}Prove Theorem \ref{thm.pie.union}, Theorem
\ref{thm.pie.nonunion}, Theorem \ref{thm.pie.jordan} and Theorem
\ref{thm.pie.bonfer}.
\end{exercise}

\subsection{Additional exercises}

This section contains some further exercises. These will not be used in the
rest of the notes, and they can be skipped at will\footnote{The same, of
course, can be said for many of the standard exercises.}. I provide solutions
to only a few of them.

\begin{exercise}
\label{exe.prop.vandermonde.consequences.f}Find a different proof of
Proposition \ref{prop.vandermonde.consequences} \textbf{(f)} that uses a
double-counting argument (i.e., counting some combinatorial objects in two
different ways, and then concluding that the results are equal).

[\textbf{Hint:} How many $\left(  x+y+1\right)  $-element subsets does the set
$\left\{  1,2,\ldots,n+1\right\}  $ have? Now, for a given $k\in\left\{
0,1,\ldots,n\right\}  $, how many $\left(  x+y+1\right)  $-element subsets
whose $\left(  x+1\right)  $-th smallest element is $k+1$ does the set
$\left\{  1,2,\ldots,n+1\right\}  $ have?]
\end{exercise}

\begin{exercise}
\label{exe.multichoose}Let $n\in\mathbb{N}$ and $k\in\mathbb{N}$ be fixed.
Show that the number of all $k$-tuples $\left(  a_{1},a_{2},\ldots
,a_{k}\right)  \in\mathbb{N}^{k}$ satisfying $a_{1}+a_{2}+\cdots+a_{k}=n$
equals $\dbinom{n+k-1}{n}$.
\end{exercise}

\begin{remark}
Exercise \ref{exe.multichoose} can be restated in terms of multisets. Namely,
let $n\in\mathbb{N}$ and $k\in\mathbb{N}$ be fixed. Also, fix a $k$-element
set $K$. Then, the number of $n$-element multisets whose elements all belong
to $K$ is $\dbinom{n+k-1}{n}$. Indeed, we can WLOG assume that $K=\left\{
1,2,\ldots,k\right\}  $ (otherwise, just relabel the elements of $K$); then,
the multisets whose elements all belong to $K$ are in bijection with the
$k$-tuples $\left(  a_{1},a_{2},\ldots,a_{k}\right)  \in\mathbb{N}^{k}$. The
bijection sends a multiset $M$ to the $k$-tuple $\left(  m_{1}\left(
M\right)  ,m_{2}\left(  M\right)  ,\ldots,m_{k}\left(  M\right)  \right)  $,
where each $m_{i}\left(  M\right)  $ is the multiplicity of the element $i$ in
$M$. The size of a multiset $M$ corresponds to the sum $a_{1}+a_{2}%
+\cdots+a_{k}$ of the entries of the resulting $k$-tuple; thus, we get a
bijection between

\begin{itemize}
\item the $n$-element multisets whose elements all belong to $K$
\end{itemize}

\noindent and

\begin{itemize}
\item the $k$-tuples $\left(  a_{1},a_{2},\ldots,a_{k}\right)  \in
\mathbb{N}^{k}$ satisfying $a_{1}+a_{2}+\cdots+a_{k}=n$.
\end{itemize}

As a consequence, Exercise \ref{exe.multichoose} shows that the number of the
former multisets is $\dbinom{n+k-1}{n}$.

Similarly, we can reinterpret the classical combinatorial interpretation of
$\dbinom{k}{n}$ (as the number of $n$-element subsets of $\left\{
1,2,\ldots,k\right\}  $) as follows: The number of all $k$-tuples $\left(
a_{1},a_{2},\ldots,a_{k}\right)  \in\left\{  0,1\right\}  ^{k}$ satisfying
$a_{1}+a_{2}+\cdots+a_{k}=n$ equals $\dbinom{k}{n}$.

See \cite[\S 13]{Galvin} and \cite[\S 1.11]{Loehr-BC} for more about multisets.
\end{remark}

\begin{exercise}
\label{exe.schmitt-iha-eq9.2}Let $n$ and $a$ be two integers with $n\geq
a\geq1$. Prove that%
\[
\sum_{k=a}^{n}\dfrac{\left(  -1\right)  ^{k}}{k}\dbinom{n-a}{k-a}%
=\dfrac{\left(  -1\right)  ^{a}}{a\dbinom{n}{a}}.
\]

\end{exercise}

\begin{exercise}
\label{exe.multichoose-app}Let $n\in\mathbb{N}$ and $k\in\mathbb{N}$. Prove
that%
\[
\sum_{u=0}^{k}\dbinom{n+u-1}{u}\dbinom{n}{k-2u}=\dbinom{n+k-1}{k}.
\]
Here, $\dbinom{a}{b}$ is defined to be $0$ when $b<0$.
\end{exercise}

Exercise \ref{exe.multichoose-app} is solved in \cite{Gri-QEDMO4P13}.

\begin{exercise}
\label{exe.bininv}Let $N\in\mathbb{N}$. The \textit{binomial transform} of a
finite sequence $\left(  f_{0},f_{1},\ldots,f_{N}\right)  \in\mathbb{Z}^{N+1}$
is defined to be the sequence $\left(  g_{0},g_{1},\ldots,g_{N}\right)  $
defined by%
\[
g_{n}=\sum_{i=0}^{n}\left(  -1\right)  ^{i}\dbinom{n}{i}f_{i}%
\ \ \ \ \ \ \ \ \ \ \text{for every }n\in\left\{  0,1,\ldots,N\right\}  .
\]

\textbf{(a)} Let $\left(  f_{0},f_{1},\ldots,f_{N}\right)  \in\mathbb{Z}%
^{N+1}$ be a finite sequence of integers. Let $\left(  g_{0},g_{1}%
,\ldots,g_{N}\right)  $ be the binomial transform of $\left(  f_{0}%
,f_{1},\ldots,f_{N}\right)  $. Show that $\left(  f_{0},f_{1},\ldots
,f_{N}\right)  $ is, in turn, the binomial transform of $\left(  g_{0}%
,g_{1},\ldots,g_{N}\right)  $.

\textbf{(b)} Find the binomial transform of the sequence $\left(
1,1,\ldots,1\right)  $.

\textbf{(c)} For any given $a\in\mathbb{N}$, find the binomial transform of
the sequence $\left(  \dbinom{0}{a},\dbinom{1}{a},\ldots,\dbinom{N}{a}\right)
$.

\textbf{(d)} For any given $q\in\mathbb{Z}$, find the binomial transform of
the sequence $\left(  q^{0},q^{1},\ldots,q^{N}\right)  $.

\textbf{(e)} Find the binomial transform of the sequence $\left(
1,0,1,0,1,0,\ldots\right)  $ (this ends with $1$ if $N$ is even, and with $0$
if $N$ is odd).

\textbf{(f)} Let $B:\mathbb{Z}^{N+1}\rightarrow\mathbb{Z}^{N+1}$ be the map
which sends every sequence $\left(  f_{0},f_{1},\ldots,f_{N}\right)
\in\mathbb{Z}^{N+1}$ to its binomial transform $\left(  g_{0},g_{1}%
,\ldots,g_{N}\right)  \in\mathbb{Z}^{N+1}$. Thus, part \textbf{(a)} of this
exercise states that $B^{2}=\operatorname*{id}$.

On the other hand, let $W:\mathbb{Z}^{N+1}\rightarrow\mathbb{Z}^{N+1}$ be the
map which sends every sequence $\left(  f_{0},f_{1},\ldots,f_{N}\right)
\in\mathbb{Z}^{N+1}$ to $\left(  \left(  -1\right)  ^{N}f_{N},\left(
-1\right)  ^{N}f_{N-1},\ldots,\left(  -1\right)  ^{N}f_{0}\right)
\in\mathbb{Z}^{N+1}$. It is rather clear that $W^{2}=\operatorname*{id}$.

Show that, furthermore, $B\circ W\circ B=W\circ B\circ W$ and $\left(  B\circ
W\right)  ^{3}=\operatorname*{id}$.
\end{exercise}

\begin{exercise}
\label{exe.binom.Hn-altsum}Let $n\in\mathbb{N}$. Prove that%
\[
\sum_{k=1}^{n}\dfrac{\left(  -1\right)  ^{k-1}}{k}\dbinom{n}{k}=\dfrac{1}%
{1}+\dfrac{1}{2}+\cdots+\dfrac{1}{n}.
\]

[\textbf{Hint:} How does the left hand side grow when $n$ is replaced by $n+1$ ?]
\end{exercise}

Exercise \ref{exe.binom.Hn-altsum} is taken from \cite[Example 3.7]{AndFen04}.

\begin{exercise}
\label{exe.binid.kurlis}Let $n\in\mathbb{N}$.

\textbf{(a)} Prove that%
\[
\sum_{k=0}^{n}\dfrac{\left(  -1\right)  ^{k}}{\dbinom{n}{k}}=2\cdot\dfrac
{n+1}{n+2}\left[  n\text{ is even}\right]  .
\]

(Here, we are using the Iverson bracket notation, as in Definition
\ref{def.iverson}; thus, $\left[  n\text{ is even}\right]  $ is $1$ if $n$ is
even and $0$ otherwise.)

\textbf{(b)} Prove that
\[
\sum_{k=0}^{n}\dfrac{1}{\dbinom{n}{k}}=\dfrac{n+1}{2^{n+1}}\sum_{k=1}%
^{n+1}\dfrac{2^{k}}{k}.
\]

[\textbf{Hint:} Show that $\dfrac{1}{\dbinom{n}{k}}=\left(  \dfrac{1}%
{\dbinom{n+1}{k}}+\dfrac{1}{\dbinom{n+1}{k+1}}\right)  \dfrac{n+1}{n+2}$ for
each $k\in\left\{  0,1,\ldots,n\right\}  $.]
\end{exercise}

Exercise \ref{exe.binid.kurlis} \textbf{(a)} is \cite[(8)]{KurLis78}. Exercise
\ref{exe.binid.kurlis} \textbf{(b)} is \cite[(9)]{KurLis78} and \cite[Example
3.9]{AndFen04} and \cite[Lemma 3.14]{AndDosS} and part of \cite[Theorem
1]{Rocket81}, and also appears with proof in
\url{https://math.stackexchange.com/a/481686/} (where it is used to show that
$\lim\limits_{n\rightarrow\infty}\sum_{k=0}^{n}\dfrac{1}{\dbinom{n}{k}}=2$).

\begin{exercise}
\label{exe.AoPS333199}For any $n\in\mathbb{N}$ and $m\in\mathbb{N}$, define a
polynomial $Z_{m,n}\in\mathbb{Z}\left[  X\right]  $ by%
\[
Z_{m,n}=\sum_{k=0}^{n}\left(  -1\right)  ^{k}\dbinom{n}{k}\left(
X^{n-k}-1\right)  ^{m}.
\]
Show that $Z_{m,n}=Z_{n,m}$ for any $n\in\mathbb{N}$ and $m\in\mathbb{N}$.
\end{exercise}

\begin{exercise}
\label{exe.AoPS262752}Let $n\in\mathbb{N}$. Prove%
\[
\sum_{k=0}^{n}\left(  -1\right)  ^{k}\dbinom{X}{k}\dbinom{X}{n-k}=%
\begin{cases}
\left(  -1\right)  ^{n/2}\dbinom{X}{n/2}, & \text{if }n\text{ is even};\\
0, & \text{if }n\text{ is odd}%
\end{cases}
\]
(an identity between polynomials in $\mathbb{Q}\left[  X\right]  $).

[\textbf{Hint:} It is enough to prove this when $X$ is replaced by a
nonnegative integer $r$ (why?). Now that you have gotten rid of polynomials,
introduce new polynomials. Namely, compute the coefficient of $X^{n}$ in
$\left(  1+X\right)  ^{r}\left(  1-X\right)  ^{r}$. Compare with the
coefficient of $X^{n}$ in $\left(  1-X^{2}\right)  ^{r}$.]
\end{exercise}

\begin{exercise}
\label{exe.vander-1/2}Let $n\in\mathbb{N}$.

\textbf{(a)} Prove that%
\[
\sum_{k=0}^{n}\dbinom{2k}{k}\dbinom{2\left(  n-k\right)  }{n-k}=4^{n}.
\]

\textbf{(b)} Prove that%
\[
\sum_{k=0}^{n}\left(  -1\right)  ^{k}\dbinom{2k}{k}\dbinom{2\left(
n-k\right)  }{n-k}=%
\begin{cases}
2^{n}\dbinom{n}{n/2}, & \text{if }n\text{ is even};\\
0, & \text{if }n\text{ is odd}%
\end{cases}
.
\]

[\textbf{Hint:} Recall Exercise \ref{exe.bin.-1/2} \textbf{(b)}.]
\end{exercise}

\begin{exercise}
\label{exe.central-binomial-even}Let $m$ be a positive integer. Prove the following:

\textbf{(a)} The binomial coefficient $\dbinom{2m}{m}$ is even.

\textbf{(b)} If $m$ is odd and satisfies $m>1$, then the binomial coefficient
$\dbinom{2m-1}{m-1}$ is even.

\textbf{(c)} If $m$ is odd and satisfies $m>1$, then $\dbinom{2m}{m}%
\equiv0\operatorname{mod}4$.
\end{exercise}

\begin{exercise}
\label{exe.supercat}For any $m\in\mathbb{N}$ and $n\in\mathbb{N}$, define a
rational number $T\left(  m,n\right)  $ by%
\[
T\left(  m,n\right)  =\dfrac{\left(  2m\right)  !\left(  2n\right)
!}{m!n!\left(  m+n\right)  !}.
\]

Prove the following facts:

\textbf{(a)} We have $4T\left(  m,n\right)  =T\left(  m+1,n\right)  +T\left(
m,n+1\right)  $ for every $m\in\mathbb{N}$ and $n\in\mathbb{N}$.

\textbf{(b)} We have $T\left(  m,n\right)  \in\mathbb{N}$ for every
$m\in\mathbb{N}$ and $n\in\mathbb{N}$.

\textbf{(c)} If $m\in\mathbb{N}$ and $n\in\mathbb{N}$ are such that $\left(
m,n\right)  \neq\left(  0,0\right)  $, then the integer $T\left(  m,n\right)
$ is even.

\textbf{(d)} If $m\in\mathbb{N}$ and $n\in\mathbb{N}$ are such that $m+n$ is
odd and $m+n>1$, then $4\mid T\left(  m,n\right)  $.

\textbf{(e)} We have $T\left(  m,0\right)  =\dbinom{2m}{m}$ for every
$m\in\mathbb{N}$.

\textbf{(f)} We have $T\left(  m,n\right)  =\dfrac{\dbinom{2m}{m}\dbinom
{2n}{n}}{\dbinom{m+n}{m}}$ for every $m\in\mathbb{N}$ and $n\in\mathbb{N}$.

\textbf{(g)} We have $T\left(  m,n\right)  =T\left(  n,m\right)  $ for every
$m\in\mathbb{N}$ and $n\in\mathbb{N}$.

\textbf{(h)} Let $m\in\mathbb{N}$ and $n\in\mathbb{N}$. Let $p=\min\left\{
m,n\right\}  $. Then,%
\[
\sum_{k=-p}^{p}\left(  -1\right)  ^{k}\dbinom{m+n}{m+k}\dbinom{m+n}%
{n+k}=\dbinom{m+n}{m}.
\]

\textbf{(i)} Let $m\in\mathbb{N}$ and $n\in\mathbb{N}$. Let $p=\min\left\{
m,n\right\}  $. Then,%
\[
T\left(  m,n\right)  =\sum_{k=-p}^{p}\left(  -1\right)  ^{k}\dbinom{2m}%
{m+k}\dbinom{2n}{n-k}.
\]

\end{exercise}

\begin{remark}
The numbers $T\left(  m,n\right)  $ introduced in Exercise \ref{exe.supercat}
are the so-called \textit{super-Catalan numbers}; much has been written about
them (e.g., \cite{Gessel92} and \cite{AleGhe14}). Exercise \ref{exe.supercat}
\textbf{(b)} suggests that these numbers count something, but no one has so
far discovered what. Exercise \ref{exe.supercat} \textbf{(i)} is a result of
von Szily (1894); see \cite[(29)]{Gessel92}. Exercise \ref{exe.supercat}
\textbf{(b)} is a result of Eug\`{e}ne Catalan (1874), and has also been posed
as \href{https://artofproblemsolving.com/community/c6h21614p139664}{Problem 3
of the International Mathematical Olympiad 1972}. Parts of Exercise
\ref{exe.supercat} are also discussed on the thread
\url{https://artofproblemsolving.com/community/c6h1553916s1_supercatalan_numbers}
.
\end{remark}

The following exercise is a variation on (\ref{eq.binom.int}):

\begin{exercise}
\label{exe.choose.a/b}Let $a$ and $b$ be two integers such that $b\neq0$. Let
$n\in\mathbb{N}$. Show that there exists some $N\in\mathbb{N}$ such that
$b^{N}\dbinom{a/b}{n}\in\mathbb{Z}$.

[\textbf{Hint:} I am not aware of a combinatorial solution to this exercise!
(I.e., I don't know what the numbers $b^{N}\dbinom{a/b}{n}$ count, even when
they are nonnegative.) All solutions that I know use some (elementary) number
theory. For the probably slickest (although unmotivated) solution, basic
modular arithmetic suffices; here is a roadmap: First, show that if $b$ and
$c$ are integers such that $c>0$, then there exists an $s\in\mathbb{Z}$ such
that $b^{c-1}\equiv sb^{c}\operatorname{mod}c$\ \ \ \ \footnotemark. Apply
this to $c=n!$ and conclude that $b^{n!}\left(  a/b-i\right)  \equiv
b^{n!}\left(  sa-i\right)  \operatorname{mod}n!$ for every $i\in\mathbb{Z}$.
Now use $\dbinom{sa}{n}\in\mathbb{Z}$.]
\end{exercise}

\footnotetext{To prove this, argue that at least two of $b^{0},b^{1}%
,\ldots,b^{c}$ are congruent modulo $c$.}

\begin{exercise}
\label{exe.ISL1975P7}\textbf{(a)} If $x$ and $y$ are two real numbers such
that $x+y=1$, and if $n\in\mathbb{N}$ and $m\in\mathbb{N}$, then prove that%
\[
x^{m+1}\sum_{k=0}^{n}\dbinom{m+k}{k}y^{k}+y^{n+1}\sum_{k=0}^{m}\dbinom{n+k}%
{k}x^{k}=1.
\]

\textbf{(b)} Let $n\in\mathbb{N}$. Prove that
\[
\sum_{k=0}^{n}\dbinom{n+k}{k}\dfrac{1}{2^{k}}=2^{n}.
\]

\end{exercise}

\begin{remark}
Exercise \ref{exe.ISL1975P7} \textbf{(a)} is Problem 7 from the IMO Shortlist
1975. It is also closely related to the Daubechies identity \cite{Zeilbe93}
(indeed, the first equality in \cite{Zeilbe93} follows by applying it to $p$,
$1-p$, $n-1$ and $n-1$ instead of $n$ and $m$). Exercise \ref{exe.ISL1975P7}
\textbf{(b)} is a fairly well-known identity for binomial coefficients (see,
e.g., \cite[(5.20)]{GKP}).
\end{remark}

\section{\label{chp.recur}Recurrent sequences}

\subsection{Basics}

Two of the most famous integer sequences defined recursively are the Fibonacci
sequence and the Lucas sequence:

\begin{itemize}
\item The
\textit{\href{https://en.wikipedia.org/wiki/Fibonacci_number}{Fibonacci
sequence}} is the sequence $\left(  f_{0},f_{1},f_{2},\ldots\right)  $ of
integers which is defined recursively by $f_{0}=0$, $f_{1}=1$, and
$f_{n}=f_{n-1}+f_{n-2}$ for all $n\geq2$. We have already introduced this
sequence in Example \ref{exa.rec-seq.fib}. Its first terms are%
\begin{align*}
f_{0}  &  =0,\ \ \ \ \ \ \ \ \ \ f_{1}=1,\ \ \ \ \ \ \ \ \ \ f_{2}%
=1,\ \ \ \ \ \ \ \ \ \ f_{3}=2,\ \ \ \ \ \ \ \ \ \ f_{4}%
=3,\ \ \ \ \ \ \ \ \ \ f_{5}=5,\\
f_{6}  &  =8,\ \ \ \ \ \ \ \ \ \ f_{7}=13,\ \ \ \ \ \ \ \ \ \ f_{8}%
=21,\ \ \ \ \ \ \ \ \ \ f_{9}=34,\ \ \ \ \ \ \ \ \ \ f_{10}=55,\\
f_{11}  &  =89,\ \ \ \ \ \ \ \ \ \ f_{12}=144,\ \ \ \ \ \ \ \ \ \ f_{13}=233.
\end{align*}
(Some authors\footnote{such as Vorobiev in his book \cite{Vorobi02}} prefer to
start the sequence at $f_{1}$ rather than $f_{0}$; of course, the recursive
definition then needs to be modified to require $f_{2}=1$ instead of $f_{0}=0$.)

\item The \textit{\href{https://en.wikipedia.org/wiki/Lucas_number}{Lucas
sequence}} is the sequence $\left(  \ell_{0},\ell_{1},\ell_{2},\ldots\right)
$ of integers which is defined recursively by $\ell_{0}=2$, $\ell_{1}=1$, and
$\ell_{n}=\ell_{n-1}+\ell_{n-2}$ for all $n\geq2$. Its first terms are%
\begin{align*}
\ell_{0}  &  =2,\ \ \ \ \ \ \ \ \ \ \ell_{1}=1,\ \ \ \ \ \ \ \ \ \ \ell
_{2}=3,\ \ \ \ \ \ \ \ \ \ \ell_{3}=4,\ \ \ \ \ \ \ \ \ \ \ell_{4}%
=7,\ \ \ \ \ \ \ \ \ \ \ell_{5}=11,\\
\ell_{6}  &  =18,\ \ \ \ \ \ \ \ \ \ \ell_{7}=29,\ \ \ \ \ \ \ \ \ \ \ell
_{8}=47,\ \ \ \ \ \ \ \ \ \ \ell_{9}=76,\ \ \ \ \ \ \ \ \ \ \ell_{10}=123,\\
\ell_{11}  &  =199,\ \ \ \ \ \ \ \ \ \ \ell_{12}=322,\ \ \ \ \ \ \ \ \ \ \ell
_{13}=521.
\end{align*}

\end{itemize}

A lot of papers and even books have been written about these two sequences,
the relations between them, and the identities that hold for their
terms.\footnote{See \url{https://oeis.org/A000045} and
\url{https://oeis.org/A000032} for an overview of their properties. The book
\cite{Vorobi02} is a readable introduction to the Fibonacci sequence, which
also surveys a lot of other mathematics (elementary number theory, continued
fractions, and even some geometry) along the way. Another introduction to the
Fibonacci sequence is \cite{CamFon07}.} One of their most striking properties
is that they can be computed explicitly, albeit using irrational numbers. In
fact, the \textit{Binet formula} says that the $n$-th Fibonacci number $f_{n}$
can be computed by%
\begin{equation}
f_{n}=\dfrac{1}{\sqrt{5}}\varphi^{n}-\dfrac{1}{\sqrt{5}}\psi^{n},
\label{eq.binet.f}%
\end{equation}
where $\varphi=\dfrac{1+\sqrt{5}}{2}$ and $\psi=\dfrac{1-\sqrt{5}}{2}$ are the
two solutions of the quadratic equation $X^{2}-X-1=0$. (The number $\varphi$
is known as the \textit{golden ratio}; the number $\psi$ can be obtained from
it by $\psi=1-\varphi=-1/\varphi$.) A similar formula, using the very same
numbers $\varphi$ and $\psi$, exists for the Lucas numbers:%
\begin{equation}
\ell_{n}=\varphi^{n}+\psi^{n}. \label{eq.binet.l}%
\end{equation}

\begin{remark}
How easy is it to compute $f_{n}$ and $\ell_{n}$ using the formulas
(\ref{eq.binet.f}) and (\ref{eq.binet.l})?

This is a nontrivial question. Indeed, if you are careless, you may find them
rather useless. For instance, if you try to compute $f_{n}$ using the formula
(\ref{eq.binet.f}) and using approximate values for the irrational numbers
$\varphi$ and $\psi$, then you might end up with a wrong value for $f_{n}$,
because the error in the approximate value for $\varphi$ propagates when you
take $\varphi$ to the $n$-th power. (And for high enough $n$, the error will
become larger than $1$, so you will not be able to get the correct value by
rounding.) The greater $n$ is, the more precise you need a value for $\varphi$
to approximate $f_{n}$ this way. Thus, approximating $\varphi$ is not a good
way to compute $f_{n}$. (Actually, the opposite is true: You can use
(\ref{eq.binet.f}) to approximate $\varphi$ by computing Fibonacci numbers.
Namely, it is easy to show that $\varphi=\lim\limits_{n\rightarrow\infty
}\dfrac{f_{n}}{f_{n-1}}$.)

A better approach to using (\ref{eq.binet.f}) is to work with the exact values
of $\varphi$ and $\psi$. To do so, you need to know how to add, subtract,
multiply and divide real numbers of the form $a+b\sqrt{5}$ with $a,b\in
\mathbb{Q}$ without ever using approximations. (Clearly, $\varphi$, $\psi$ and
$\sqrt{5}$ all have this form.) There are rules for this, which are simple to
check:%
\begin{align*}
\left(  a+b\sqrt{5}\right)  +\left(  c+d\sqrt{5}\right)   &  =\left(
a+c\right)  +\left(  b+d\right)  \sqrt{5};\\
\left(  a+b\sqrt{5}\right)  -\left(  c+d\sqrt{5}\right)   &  =\left(
a-c\right)  +\left(  b-d\right)  \sqrt{5};\\
\left(  a+b\sqrt{5}\right)  \cdot\left(  c+d\sqrt{5}\right)   &  =\left(
ac+5bd\right)  +\left(  bc+ad\right)  \sqrt{5};\\
\dfrac{a+b\sqrt{5}}{c+d\sqrt{5}}  &  =\dfrac{\left(  ac-5bd\right)  +\left(
bc-ad\right)  \sqrt{5}}{c^{2}-5d^{2}}\ \ \ \ \ \ \ \ \ \ \text{for }\left(
c,d\right)  \neq\left(  0,0\right)  .
\end{align*}
(The last rule is an instance of \textquotedblleft rationalizing the
denominator\textquotedblright.) These rules give you a way to exactly compute
things like $\varphi^{n}$, $\dfrac{1}{\sqrt{5}}\varphi^{n}$, $\psi^{n}$ and
$\dfrac{1}{\sqrt{5}}\psi^{n}$, and thus also $f_{n}$ and $\ell_{n}$. If you
use
\href{https://en.wikipedia.org/wiki/Exponentiation_by_squaring}{exponentiation
by squaring} to compute $n$-th powers, this actually becomes a fast algorithm
(a lot faster than just computing $f_{n}$ and $\ell_{n}$ using the
recurrence). So, yes, (\ref{eq.binet.f}) and (\ref{eq.binet.l}) are useful.
\end{remark}

We shall now study a generalization of both the Fibonacci and the Lucas
sequences, and generalize (\ref{eq.binet.f}) and (\ref{eq.binet.l}) to a
broader class of sequences.

\begin{definition}
\label{def.abrec} If $a$ and $b$ are two complex numbers, then a sequence
$\left(  x_{0},x_{1},x_{2},\ldots\right)  $ of complex numbers will be called
$\left(  a,b\right)  $\textit{-recurrent} if every $n\geq2$ satisfies%
\[
x_{n}=ax_{n-1}+bx_{n-2}.
\]

\end{definition}

So, the Fibonacci sequence and the Lucas sequence are $\left(  1,1\right)
$-recurrent. An $\left(  a,b\right)  $-recurrent sequence $\left(  x_{0}%
,x_{1},x_{2},\ldots\right)  $ is fully determined by the four values $a$, $b$,
$x_{0}$ and $x_{1}$, and can be constructed for any choice of these four
values. Here are some further examples of $\left(  a,b\right)  $-recurrent sequences:

\begin{itemize}
\item The sequence $\left(  x_{0},x_{1},x_{2},\ldots\right)  $ in Theorem
\ref{thm.rec-seq.fibx} is $\left(  a,b\right)  $-recurrent (by its very definition).

\item A sequence $\left(  x_{0},x_{1},x_{2},\ldots\right)  $ is $\left(
2,-1\right)  $-recurrent if and only if every $n\geq2$ satisfies
$x_{n}=2x_{n-1}-x_{n-2}$. In other words, a sequence $\left(  x_{0}%
,x_{1},x_{2},\ldots\right)  $ is $\left(  2,-1\right)  $-recurrent if and only
if every $n\geq2$ satisfies $x_{n}-x_{n-1}=x_{n-1}-x_{n-2}$. In other words, a
sequence $\left(  x_{0},x_{1},x_{2},\ldots\right)  $ is $\left(  2,-1\right)
$-recurrent if and only if $x_{1}-x_{0}=x_{2}-x_{1}=x_{3}-x_{2}=\cdots$. In
other words, the $\left(  2,-1\right)  $-recurrent sequences are precisely the
arithmetic progressions.

\item Geometric progressions are also $\left(  a,b\right)  $-recurrent for
appropriate $a$ and $b$. Namely, any geometric progression $\left(
u,uq,uq^{2},uq^{3},\ldots\right)  $ is $\left(  q,0\right)  $-recurrent, since
every $n\geq2$ satisfies $uq^{n}=q\cdot uq^{n-1}+0\cdot uq^{n-2}$. However,
not every $\left(  q,0\right)  $-recurrent sequence $\left(  x_{0},x_{1}%
,x_{2},\ldots\right)  $ is a geometric progression (since the condition
$x_{n}=qx_{n-1}+0x_{n-2}$ for all $n\geq2$ says nothing about $x_{0}$, and
thus $x_{0}$ can be arbitrary).

\item A sequence $\left(  x_{0},x_{1},x_{2},\ldots\right)  $ is $\left(
0,1\right)  $-recurrent if and only if every $n\geq2$ satisfies $x_{n}%
=x_{n-2}$. In other words, a sequence $\left(  x_{0},x_{1},x_{2}%
,\ldots\right)  $ is $\left(  0,1\right)  $-recurrent if and only if it has
the form $\left(  u,v,u,v,u,v,\ldots\right)  $ for two complex numbers $u$ and
$v$.

\item A sequence $\left(  x_{0},x_{1},x_{2},\ldots\right)  $ is $\left(
1,0\right)  $-recurrent if and only if every $n\geq2$ satisfies $x_{n}%
=x_{n-1}$. In other words, a sequence $\left(  x_{0},x_{1},x_{2}%
,\ldots\right)  $ is $\left(  1,0\right)  $-recurrent if and only if it has
the form $\left(  u,v,v,v,v,\ldots\right)  $ for two complex numbers $u$ and
$v$. Notice that $u$ is not required to be equal to $v$, because we never
claimed that $x_{n}=x_{n-1}$ holds for $n=1$.

\item A sequence $\left(  x_{0},x_{1},x_{2},\ldots\right)  $ is $\left(
1,-1\right)  $-recurrent if and only if every $n\geq2$ satisfies
$x_{n}=x_{n-1}-x_{n-2}$. Curiously, it turns out that every such sequence is
$6$-periodic (i.e., it satisfies $x_{n+6}=x_{n}$ for every $n\in\mathbb{N}$),
because every $n\in\mathbb{N}$ satisfies%
\begin{align*}
x_{n+6}  &  =\underbrace{x_{n+5}}_{=x_{n+4}-x_{n+3}}-x_{n+4}=\left(
x_{n+4}-x_{n+3}\right)  -x_{n+4}=-\underbrace{x_{n+3}}_{=x_{n+2}-x_{n+1}}\\
&  =-\left(  \underbrace{x_{n+2}}_{=x_{n+1}-x_{n}}-x_{n+1}\right)  =-\left(
x_{n+1}-x_{n}-x_{n+1}\right)  =x_{n}.
\end{align*}
More precisely, a sequence $\left(  x_{0},x_{1},x_{2},\ldots\right)  $ is
$\left(  1,-1\right)  $-recurrent if and only if it has the form $\left(
u,v,v-u,-u,-v,u-v,\ldots\right)  $ (where the \textquotedblleft$\ldots
$\textquotedblright\ stands for \textquotedblleft repeat the preceding $6$
values over and over\textquotedblright\ here) for two complex numbers $u$ and
$v$.

\item The above three examples notwithstanding, most $\left(  a,b\right)
$-recurrent sequences of course are not periodic. However, here is another
example which provides a great supply of non-periodic $\left(  a,b\right)
$-recurrent sequences and, at the same time, explains why we get so many
periodic ones: If $\alpha$ is any angle, then the sequences%
\begin{align*}
&  \left(  \sin\left(  0\alpha\right)  ,\sin\left(  1\alpha\right)
,\sin\left(  2\alpha\right)  ,\ldots\right)  \ \ \ \ \ \ \ \ \ \ \text{and}\\
&  \left(  \cos\left(  0\alpha\right)  ,\cos\left(  1\alpha\right)
,\cos\left(  2\alpha\right)  ,\ldots\right)
\end{align*}
are $\left(  2\cos\alpha,-1\right)  $-recurrent. More generally, if $\alpha$
and $\beta$ are two angles, then the sequence%
\[
\left(  \sin\left(  \beta+0\alpha\right)  ,\sin\left(  \beta+1\alpha\right)
,\sin\left(  \beta+2\alpha\right)  ,\ldots\right)
\]
is $\left(  2\cos\alpha,-1\right)  $-recurrent\footnote{\textit{Proof.} Let
$\alpha$ and $\beta$ be two angles. We need to show that the sequence $\left(
\sin\left(  \beta+0\alpha\right)  ,\sin\left(  \beta+1\alpha\right)
,\sin\left(  \beta+2\alpha\right)  ,\ldots\right)  $ is $\left(  2\cos
\alpha,-1\right)  $-recurrent. In other words, we need to prove that%
\[
\sin\left(  \beta+n\alpha\right)  =2\cos\alpha\sin\left(  \beta+\left(
n-1\right)  \alpha\right)  +\left(  -1\right)  \sin\left(  \beta+\left(
n-2\right)  \alpha\right)
\]
for every $n\geq2$. So fix $n\geq2$.
\par
One of the well-known trigonometric identities states that $\sin x+\sin
y=2\sin\dfrac{x+y}{2}\cos\dfrac{x-y}{2}$ for any two angles $x$ and $y$.
Applying this to $x=\beta+n\alpha$ and $y=\beta+\left(  n-2\right)  \alpha$,
we obtain%
\begin{align*}
\sin\left(  \beta+n\alpha\right)  +\sin\left(  \beta+\left(  n-2\right)
\alpha\right)   &  =2\sin\underbrace{\dfrac{\left(  \beta+n\alpha\right)
+\left(  \beta+\left(  n-2\right)  \alpha\right)  }{2}}_{=\beta+\left(
n-1\right)  \alpha}\cos\underbrace{\dfrac{\left(  \beta+n\alpha\right)
-\left(  \beta+\left(  n-2\right)  \alpha\right)  }{2}}_{=\alpha}\\
&  =2\sin\left(  \beta+\left(  n-1\right)  \alpha\right)  \cos\alpha
=2\cos\alpha\sin\left(  \beta+\left(  n-1\right)  \alpha\right)  .
\end{align*}
Hence,
\begin{align*}
\sin\left(  \beta+n\alpha\right)   &  =2\cos\alpha\sin\left(  \beta+\left(
n-1\right)  \alpha\right)  -\sin\left(  \beta+\left(  n-2\right)
\alpha\right) \\
&  =2\cos\alpha\sin\left(  \beta+\left(  n-1\right)  \alpha\right)  +\left(
-1\right)  \sin\left(  \beta+\left(  n-2\right)  \alpha\right)  ,
\end{align*}
qed.}. When $\alpha\in2\pi\mathbb{Q}$ (that is, $\alpha=2\pi r$ for some
$r\in\mathbb{Q}$), this sequence is periodic.
\end{itemize}

\subsection{Explicit formulas (\`{a} la Binet)}

Now, we can get an explicit formula (similar to (\ref{eq.binet.f}) and
(\ref{eq.binet.l})) for every term of an $\left(  a,b\right)  $-recurrent
sequence (in terms of $a$, $b$, $x_{0}$ and $x_{1}$) in the case when
$a^{2}+4b\neq0$. Here is how this works:

\begin{remark}
\label{rmk.binet}Let $a$ and $b$ be complex numbers such that $a^{2}+4b\neq0$.
Let $\left(  x_{0},x_{1},x_{2},\ldots\right)  $ be an $\left(  a,b\right)
$-recurrent sequence. We want to construct an explicit formula for each
$x_{n}$ in terms of $x_{0}$, $x_{1}$, $a$ and $b$.

To do so, we let $q_{+}$ and $q_{-}$ be the two solutions of the quadratic
equation $X^{2}-aX-b=0$, namely%
\[
q_{+}=\dfrac{a+\sqrt{a^{2}+4b}}{2}\ \ \ \ \ \ \ \ \ \ \text{and}%
\ \ \ \ \ \ \ \ \ \ q_{-}=\dfrac{a-\sqrt{a^{2}+4b}}{2}.
\]
We notice that $q_{+}\neq q_{-}$ (since $a^{2}+4b\neq0$). It is easy to see
that the sequences $\left(  1,q_{+},q_{+}^{2},q_{+}^{3},\ldots\right)  $ and
$\left(  1,q_{-},q_{-}^{2},q_{-}^{3},\ldots\right)  $ are $\left(  a,b\right)
$-recurrent. As a consequence, for any two complex numbers $\lambda_{+}$ and
$\lambda_{-}$, the sequence%
\[
\left(  \lambda_{+}+\lambda_{-},\lambda_{+}q_{+}+\lambda_{-}q_{-},\lambda
_{+}q_{+}^{2}+\lambda_{-}q_{-}^{2},\ldots\right)
\]
(the $n$-th term of this sequence, with $n$ starting at $0$, is $\lambda
_{+}q_{+}^{n}+\lambda_{-}q_{-}^{n}$) must also be $\left(  a,b\right)
$-recurrent (check this!). We denote this sequence by $L_{\lambda_{+}%
,\lambda_{-}}$.

We now need to find two complex numbers $\lambda_{+}$ and $\lambda_{-}$ such
that this sequence $L_{\lambda_{+},\lambda_{-}}$ is our sequence $\left(
x_{0},x_{1},x_{2},\ldots\right)  $. In order to do so, we only need to ensure
that $\lambda_{+}+\lambda_{-}=x_{0}$ and $\lambda_{+}q_{+}+\lambda_{-}%
q_{-}=x_{1}$ (because once this holds, it will follow that the sequences
$L_{\lambda_{+},\lambda_{-}}$ and $\left(  x_{0},x_{1},x_{2},\ldots\right)  $
have the same first two terms; and this will yield that these two sequences
are identical, because two $\left(  a,b\right)  $-recurrent sequences with the
same first two terms must be identical). That is, we need to solve the system
of linear equations%
\[
\left\{
\begin{array}
[c]{c}%
\lambda_{+}+\lambda_{-}=x_{0};\\
\lambda_{+}q_{+}+\lambda_{-}q_{-}=x_{1}%
\end{array}
\right.  \ \ \ \ \ \ \ \ \ \ \text{in the unknowns }\lambda_{+}\text{ and
}\lambda_{-}.
\]
Thanks to $q_{+}\neq q_{-}$, this system has a unique solution:%
\[
\lambda_{+}=\dfrac{x_{1}-q_{-}x_{0}}{q_{+}-q_{-}};\ \ \ \ \ \ \ \ \ \ \lambda
_{-}=\dfrac{q_{+}x_{0}-x_{1}}{q_{+}-q_{-}}.
\]
Thus, if we set $\left(  \lambda_{+},\lambda_{-}\right)  $ to be this
solution, then $\left(  x_{0},x_{1},x_{2},\ldots\right)  =L_{\lambda
_{+},\lambda_{-}}$, so that%
\begin{equation}
x_{n}=\lambda_{+}q_{+}^{n}+\lambda_{-}q_{-}^{n}
\label{rmk.recursive.binet-general}%
\end{equation}
for every nonnegative integer $n$. This is an explicit formula, at least if
the square roots do not disturb you. When $x_{0}=0$ and $x_{1}=a=b=1$, you get
the famous Binet formula (\ref{eq.binet.f}) for the Fibonacci sequence.
\end{remark}

In the next exercise you will see what happens if the $a^{2}+4b\neq0$
condition does not hold.

\begin{exercise}
\label{exe.ps2.2.1}Let $a$ and $b$ be complex numbers such that $a^{2}+4b=0$.
Consider an $\left(  a,b\right)  $-recurrent sequence $\left(  x_{0}%
,x_{1},x_{2},\ldots\right)  $. Find an explicit formula for each $x_{n}$ in
terms of $x_{0}$, $x_{1}$, $a$ and $b$.

[\textbf{Note:} The polynomial $X^{2}-aX-b$ has a double root here. Unlike the
case of two distinct roots studied above, you won't see any radicals here. The
explicit formula really deserves the name \textquotedblleft
explicit\textquotedblright.]
\end{exercise}

Remark \ref{rmk.binet} and Exercise \ref{exe.ps2.2.1}, combined, solve the
problem of finding an explicit formula for any term of an $\left(  a,b\right)
$-recurrent sequence when $a$ and $b$ are complex numbers, at least if you
don't mind having square roots in your formula. Similar tactics can be used to
find explicit forms for the more general case of sequences satisfying
\textquotedblleft homogeneous linear recurrences with constant
coefficients\textquotedblright\footnote{These are sequences $\left(
x_{0},x_{1},x_{2},\ldots\right)  $ which satisfy%
\[
\left(  x_{n}=c_{1}x_{n-1}+c_{2}x_{n-2}+\cdots+c_{k}x_{n-k}%
\ \ \ \ \ \ \ \ \ \ \text{for all }n\geq k\right)
\]
for a fixed $k\in\mathbb{N}$ and a fixed $k$-tuple $\left(  c_{1},c_{2}%
,\ldots,c_{k}\right)  $ of complex numbers. When $k=2$, these are the $\left(
c_{1},c_{2}\right)  $-recurrent sequences.}, although instead of square roots
you will now need roots of higher-degree polynomials. (See \cite[\S 22.3.2
(\textquotedblleft Solving Homogeneous Linear Recurrences\textquotedblright%
)]{LeLeMe16} for an outline of this; see also \cite[Topic \textquotedblleft
Linear Recurrences\textquotedblright]{Hefferon} for a linear-algebraic introduction.)

\subsection{Further results}

Here are some more exercises from the theory of recurrent sequences. I am not
going particularly deep here, but we may encounter generalizations later.

First, an example: If we \textquotedblleft split\textquotedblright\ the
Fibonacci sequence
\[
\left(  f_{0},f_{1},f_{2},\ldots\right)  =\left(  0,1,1,2,3,5,8,\ldots\right)
\]
into two subsequences
\[
\left(  f_{0},f_{2},f_{4},\ldots\right)  =\left(  0,1,3,8,21,\ldots\right)
\qquad\text{and}\qquad\left(  f_{1},f_{3},f_{5},\ldots\right)  =\left(
1,2,5,13,\ldots\right)
\]
(each of which contains every other Fibonacci number), then it turns out that
each of these two subsequences is $\left(  3,-1\right)  $%
-recurrent\footnote{In other words, we have $f_{2n}=3f_{2\left(  n-1\right)
}+\left(  -1\right)  f_{2\left(  n-2\right)  }$ and $f_{2n+1}=3f_{2\left(
n-1\right)  +1}+\left(  -1\right)  f_{2\left(  n-2\right)  +1}$ for every
$n\geq2$.}. This is rather easy to prove, but one can always ask for
generalizations: What happens if we start with an arbitrary $\left(
a,b\right)  $-recurrent sequence, instead of the Fibonacci numbers? What
happens if we split it into three, four or more subsequences? The answer is
rather nice:

\begin{exercise}
\label{exe.ps2.2.2}Let $a$ and $b$ be complex numbers. Let $\left(
x_{0},x_{1},x_{2},\ldots\right)  $ be an $\left(  a,b\right)  $-recurrent sequence.

\textbf{(a)} Prove that the sequences $\left(  x_{0},x_{2},x_{4}%
,\ldots\right)  $ and $\left(  x_{1},x_{3},x_{5},\ldots\right)  $ are $\left(
c,d\right)  $-recurrent for some complex numbers $c$ and $d$. Find these $c$
and $d$.

\textbf{(b)} Prove that the sequences $\left(  x_{0},x_{3},x_{6}%
,\ldots\right)  $, $\left(  x_{1},x_{4},x_{7},\ldots\right)  $ and $\left(
x_{2},x_{5},x_{8},\ldots\right)  $ are $\left(  c,d\right)  $-recurrent for
some (other) complex numbers $c$ and $d$.

\textbf{(c)} For every nonnegative integers $N$ and $K$, prove that the
sequence $\left(  x_{K},x_{N+K},x_{2N+K},x_{3N+K},\ldots\right)  $ is $\left(
c,d\right)  $-recurrent for some complex numbers $c$ and $d$ which depend only
on $N$, $a$ and $b$ (but not on $K$ or $x_{0}$ or $x_{1}$).
\end{exercise}

The next exercise gives a combinatorial interpretation of the Fibonacci numbers:

\begin{exercise}
\label{exe.ps2.2.3}Recall that the Fibonacci numbers $f_{0},f_{1},f_{2}%
,\ldots$ are defined recursively by $f_{0}=0$, $f_{1}=1$ and $f_{n}%
=f_{n-1}+f_{n-2}$ for all $n\geq2$. For every positive integer $n$, show that
$f_{n}$ is the number of subsets $I$ of $\left\{  1,2,\ldots,n-2\right\}  $
such that no two elements of $I$ are consecutive (i.e., there exists no
$i\in\mathbb{Z}$ such that both $i$ and $i+1$ belong to $I$). For instance,
for $n=5$, these subsets are $\varnothing$, $\left\{  1\right\}  $, $\left\{
2\right\}  $, $\left\{  3\right\}  $ and $\left\{  1,3\right\}  $.
\end{exercise}

Notice that $\left\{  1,2,\ldots,-1\right\}  $ is to be understood as the
empty set (since there are no integers $x$ satisfying $1\leq x\leq-1$). (So
Exercise \ref{exe.ps2.2.3}, applied to $n=1$, says that $f_{1}$ is the number
of subsets $I$ of the empty set such that no two elements of $I$ are
consecutive. This is correct, because the empty set has only one subset, which
of course is empty and thus has no consecutive elements; and the Fibonacci
number $f_{1}$ is precisely $1$.)

\begin{remark}
\label{rmk.fib.dominos}Exercise \ref{exe.ps2.2.3} is equivalent to another
known combinatorial interpretation of the Fibonacci numbers.

Namely, let $n$ be a positive integer. Consider a rectangular table of
dimensions $2\times\left(  n-1\right)  $ (that is, with $2$ rows and $n-1$
columns). How many ways are there to subdivide this table into dominos? (A
\textit{domino} means a set of two adjacent boxes.)

For $n=5$, there are $5$ ways:%
\begin{align*}
&
\begin{tabular}{ | c | c | c | c | }
\hline\multicolumn{1}{| c}{\phantom{a}} & \phantom{a} & \multicolumn{1}%
{c}{\phantom{a}} & \phantom{a} \\ \hline\multicolumn{1}{| c}{\phantom{a}}
& \phantom{a} & \multicolumn{1}{c}{\phantom{a}} & \phantom{a} \\ \hline
\end{tabular}%
\ ,\ \ \ \ \ \ \ \ \ \
\begin{tabular}{ | c | c | c | c | }
\hline\phantom{a} & \phantom{a} & \multicolumn{1}{c}{\phantom{a}}
& \phantom{a} \\ \cline{3-4}
\phantom{a} & \phantom{a} & \multicolumn{1}{c}{\phantom{a}} & \phantom{a}
\\ \hline\end{tabular}%
\ ,\ \ \ \ \ \ \ \ \ \
\begin{tabular}{ | c | c | c | c | }
\hline\multicolumn{1}{| c}{\phantom{a}} & \phantom{a} & \phantom{a}
& \phantom{a} \\ \cline{1-2}
\multicolumn{1}{| c}{\phantom{a}} & \phantom{a} & \phantom{a} & \phantom{a}
\\ \hline\end{tabular}%
\ ,\\
&
\begin{tabular}{ | c | c | c | c | }
\hline\phantom{a} & \multicolumn{1}{ c}{\phantom{a}} & \phantom{a}
& \phantom{a} \\ \cline{2-3}
\phantom{a} & \multicolumn{1}{ c}{\phantom{a}} & \phantom{a} & \phantom{a}
\\ \hline\end{tabular}%
\ ,\ \ \ \ \ \ \ \ \ \
\begin{tabular}{ | c | c | c | c | }
\hline\phantom{a} & \phantom{a} & \phantom{a} & \phantom{a} \\
\phantom{a} & \phantom{a} & \phantom{a} & \phantom{a} \\ \hline\end{tabular}%
\ .
\end{align*}
In the general case, there are $f_{n}$ ways. Why?

As promised, this result is equivalent to Exercise \ref{exe.ps2.2.3}. Let us
see why. Let $P$ be a way to subdivide the table into dominos. We say that a
\textit{horizontal domino} is a domino which consists of two adjacent boxes in
the same row; similarly, we define a vertical domino. It is easy to see that
(in the subdivision $P$) each column of the table is covered either by a
single vertical domino, or by two horizontal dominos (in which case either
both of them \textquotedblleft begin\textquotedblright\ in this column, or
both of them \textquotedblleft end\textquotedblright\ in this column). Let
$J\left(  P\right)  $ be the set of all $i\in\left\{  1,2,\ldots,n-1\right\}
$ such that the $i$-th column of the table is covered by two horizontal
dominos, both of which \textquotedblleft begin\textquotedblright\ in this
column. For instance,%
\begin{align*}
J\left(  \
\begin{tabular}{ | c | c | c | c | }
\hline\multicolumn{1}{| c}{\phantom{a}} & \phantom{a} & \multicolumn{1}%
{c}{\phantom{a}} & \phantom{a} \\ \hline\multicolumn{1}{| c}{\phantom{a}}
& \phantom{a} & \multicolumn{1}{c}{\phantom{a}} & \phantom{a} \\ \hline
\end{tabular}%
\ \right)   &  =\left\{  1,3\right\}  ;\\
J\left(  \
\begin{tabular}{ | c | c | c | c | }
\hline\phantom{a} & \phantom{a} & \multicolumn{1}{c}{\phantom{a}}
& \phantom{a} \\ \cline{3-4}
\phantom{a} & \phantom{a} & \multicolumn{1}{c}{\phantom{a}} & \phantom{a}
\\ \hline\end{tabular}%
\ \right)   &  =\left\{  3\right\}  ;\\
J\left(  \
\begin{tabular}{ | c | c | c | c | }
\hline\multicolumn{1}{| c}{\phantom{a}} & \phantom{a} & \phantom{a}
& \phantom{a} \\ \cline{1-2}
\multicolumn{1}{| c}{\phantom{a}} & \phantom{a} & \phantom{a} & \phantom{a}
\\ \hline\end{tabular}%
\ \right)   &  =\left\{  1\right\}  ;\\
J\left(  \
\begin{tabular}{ | c | c | c | c | }
\hline\phantom{a} & \multicolumn{1}{ c}{\phantom{a}} & \phantom{a}
& \phantom{a} \\ \cline{2-3}
\phantom{a} & \multicolumn{1}{ c}{\phantom{a}} & \phantom{a} & \phantom{a}
\\ \hline\end{tabular}%
\ \right)   &  =\left\{  2\right\}  ;\\
J\left(  \
\begin{tabular}{ | c | c | c | c | }
\hline\phantom{a} & \phantom{a} & \phantom{a} & \phantom{a} \\
\phantom{a} & \phantom{a} & \phantom{a} & \phantom{a} \\ \hline\end{tabular}%
\ \right)   &  =\varnothing.
\end{align*}
It is easy to see that the set $J\left(  P\right)  $ is a subset of $\left\{
1,2,\ldots,n-2\right\}  $ containing no two consecutive integers. Moreover,
this set $J\left(  P\right)  $ uniquely determines $P$, and for every subset
$I$ of $\left\{  1,2,\ldots,n-2\right\}  $ containing no two consecutive
integers, there exists some way $P$ to subdivide the table into dominos such
that $J\left(  P\right)  =I$.

Hence, the number of all ways to subdivide the table into dominos equals the
number of all subsets $I$ of $\left\{  1,2,\ldots,n-2\right\}  $ containing no
two consecutive integers. Exercise \ref{exe.ps2.2.3} says that this latter
number is $f_{n}$; therefore, so is the former number.

(I have made this remark because I found it instructive. If you merely want a
proof that the number of all ways to subdivide the table into dominos equals
$f_{n}$, then I guess it is easier to just prove it by induction without
taking the detour through Exercise \ref{exe.ps2.2.3}. This proof is sketched
in \cite[\S 7.1]{GKP}, followed by an informal yet insightful discussion of
\textquotedblleft infinite sums of dominos\textquotedblright\ and various
related ideas.)
\end{remark}

Either Exercise \ref{exe.ps2.2.3} or Remark \ref{rmk.fib.dominos} can be used
to prove properties of Fibonacci numbers in a combinatorial way; see
\cite{BenQui-fib} for some examples of such proofs.

Here is another formula for certain recursive sequences, coming out of a
recent paper on cluster algebras\footnote{Specifically, Exercise
\ref{exe.ps2.2.S} is part of \cite[Definition 1]{LS2}, but I have reindexed
the sequence and fixed the missing upper bound in the sum.}:

\Needspace{4cm}

\begin{exercise}
\label{exe.ps2.2.S}Let $r\in\mathbb{Z}$. Define a sequence $\left(
c_{0},c_{1},c_{2},\ldots\right)  $ of integers recursively by $c_{0}=0$,
$c_{1}=1$ and $c_{n}=rc_{n-1}-c_{n-2}$ for all $n\geq2$. Show that%
\begin{equation}
c_{n}=\sum_{i=0}^{n-1}\left(  -1\right)  ^{i}\dbinom{n-1-i}{i}r^{n-1-2i}
\label{eq.exe.2.S}%
\end{equation}
for every $n\in\mathbb{N}$. Here, we use the following convention: Any
expression of the form $a\cdot b$, where $a$ is $0$, has to be interpreted as
$0$, even if $b$ is undefined.\footnotemark
\end{exercise}

\footnotetext{The purpose of this convention is to make sure that the right
hand side of (\ref{eq.exe.2.S}) is well-defined, even though the expression
$r^{n-1-2i}$ that appears in it might be undefined (it will be undefined when
$r=0$ and $n-1-2i<0$).
\par
Of course, the downside of this convention is that we might not have $a\cdot
b=b\cdot a$ (because $a\cdot b$ might be well-defined while $b\cdot a$ is not,
or vice versa).}

\subsection{Additional exercises}

This section contains some further exercises. As the earlier \textquotedblleft
additional exercises\textquotedblright, these will not be relied on in the
rest of this text, and solutions will not be provided.

\begin{exercise}
\label{exe.rec.qn-rn}Let $q$ and $r$ be two complex numbers. Prove that the
sequence $\left(  q^{0}-r^{0},q^{1}-r^{1},q^{2}-r^{2},\ldots\right)  $ is
$\left(  a,b\right)  $-recurrent for two appropriately chosen $a$ and $b$.
Find these $a$ and $b$.
\end{exercise}

\begin{exercise}
\label{exe.fib.floor}Let $\varphi$ be the golden ratio (i.e., the real number
$\dfrac{1+\sqrt{5}}{2}$). Let $\left(  f_{0},f_{1},f_{2},\ldots\right)  $ be
the Fibonacci sequence.

\textbf{(a)} Show that $f_{n+1}-\varphi f_{n}=\dfrac{1}{\sqrt{5}}\psi^{n}$ for
every $n\in\mathbb{N}$, where $\psi=\dfrac{1-\sqrt{5}}{2}$. (Notice that
$\psi=\dfrac{1-\sqrt{5}}{2}\approx-0.618$ lies between $-1$ and $0$, and thus
the powers $\psi^{n}$ converge to $0$ as $n\rightarrow\infty$. So
$f_{n+1}-\varphi f_{n}\rightarrow0$ as $n\rightarrow\infty$, and consequently
$\dfrac{f_{n+1}}{f_{n}}\rightarrow\varphi$ as well.)

\textbf{(b)} Show that%
\[
f_{n}=\operatorname*{round}\left(  \dfrac{1}{\sqrt{5}}\varphi^{n}\right)
\ \ \ \ \ \ \ \ \ \ \text{for every }n\in\mathbb{N}.
\]
Here, if $x$ is a real number, then $\operatorname*{round}x$ denotes the
integer closest to $x$ (where, in case of a tie, we take the higher of the two candidates\footnotemark).
\end{exercise}

\footnotetext{This does not really matter in our situation, because $\dfrac
{1}{\sqrt{5}}\varphi^{n}$ will never be a half-integer.}

\begin{exercise}
\label{exe.fib.zeckendorf}Let $\left(  f_{0},f_{1},f_{2},\ldots\right)  $ be
the Fibonacci sequence. A set $I$ of integers is said to be \textit{lacunar}
if no two elements of $I$ are consecutive (i.e., there exists no $i\in I$ such
that $i+1\in I$). Show that, for every $n\in\mathbb{N}$, there exists a unique
lacunar subset $S$ of $\left\{  2,3,4,\ldots\right\}  $ such that
$n=\sum_{s\in S}f_{s}$.

(For example, if $n=17$, then $S=\left\{  2,4,7\right\}  $, because
$17=1+3+13=f_{2}+f_{4}+f_{7}$.)
\end{exercise}

\begin{remark}
The representation of $n$ in the form $n=\sum_{s\in S}f_{s}$ in Exercise
\ref{exe.fib.zeckendorf} is known as the
\textit{\href{https://en.wikipedia.org/wiki/Zeckendorf's_theorem}{\textit{Zeckendorf
representation}}} of $n$. It has a number of interesting properties and trivia
related to it; for example, there is
\href{https://www.encyclopediaofmath.org/index.php/Zeckendorf_representation}{a
rule of thumb for converting miles into kilometers that uses it}. It can also
be used to define a curious \textquotedblleft Fibonacci
multiplication\textquotedblright\ operation on nonnegative integers
\cite{Knuth-fib}.
\end{remark}

\begin{exercise}
\label{exe.fib.zeckendids}Let $\left(  f_{0},f_{1},f_{2},\ldots\right)  $ be
the Fibonacci sequence.

\textbf{(a)} Prove the identities%
\begin{align*}
1f_{n}  &  =f_{n}\ \ \ \ \ \ \ \ \ \ \text{for all }n\geq0;\\
2f_{n}  &  =f_{n-2}+f_{n+1}\ \ \ \ \ \ \ \ \ \ \text{for all }n\geq2;\\
3f_{n}  &  =f_{n-2}+f_{n+2}\ \ \ \ \ \ \ \ \ \ \text{for all }n\geq2;\\
4f_{n}  &  =f_{n-2}+f_{n}+f_{n+2}\ \ \ \ \ \ \ \ \ \ \text{for all }n\geq2.
\end{align*}

\textbf{(b)} Notice that the right hand sides of these identities have a
specific form: they are sums of $f_{n+t}$ for $t$ ranging over a lacunar
subset of $\mathbb{Z}$. (See Exercise \ref{exe.fib.zeckendorf} for the
definition of \textquotedblleft lacunar\textquotedblright.) Try to find
similar identities for $5f_{n}$ and $6f_{n}$.

\textbf{(c)} Prove that such identities exist in general. More precisely,
prove the following: Let $T$ be a finite set, and $a_{t}$ be an integer for
every $t\in T$. Then, there exists a unique finite lacunar subset $S$ of
$\mathbb{Z}$ such that
\begin{align*}
\sum\limits_{t\in T}f_{n+a_{t}}  &  =\sum\limits_{s\in S}f_{n+s}%
\ \ \ \ \ \ \ \ \ \ \text{for every }n\in\mathbb{Z}\text{ which}\\
&  \ \ \ \ \ \ \ \ \ \ \ \ \ \ \ \ \ \ \ \ \text{ satisfies }n\geq\max\left(
\left\{  -a_{t}\mid t\in T\right\}  \cup\left\{  -s\mid s\in S\right\}
\right)  .
\end{align*}
(The condition $n\geq\max\left(  \left\{  -a_{t}\mid t\in T\right\}
\cup\left\{  -s\mid s\in S\right\}  \right)  $ merely ensures that all the
$f_{n+a_{t}}$ and $f_{n+s}$ are well-defined.)
\end{exercise}

\begin{remark}
Exercise \ref{exe.fib.zeckendids} \textbf{(c)} is \cite[Theorem 1.4]%
{Gri-zeck}. It is also a consequence of \cite[Lemma 6.2]{CamFon07} (applied to
$k=2$). I'd be delighted to see other proofs!

Similarly I am highly interested in analogues of Exercises
\ref{exe.fib.zeckendorf} and \ref{exe.fib.zeckendids} for other $\left(
a,b\right)  $-recurrent sequences (e.g., Lucas numbers).
\end{remark}

\begin{exercise}
\label{exe.rec.fibonomial}\textbf{(a)} Let $\left(  f_{0},f_{1},f_{2}%
,\ldots\right)  $ be the Fibonacci sequence. For every $n\in\mathbb{N}$ and
$k\in\mathbb{N}$ satisfying $0\leq k\leq n$, define a rational number
$\dbinom{n}{k}_{F}$ by%
\[
\dbinom{n}{k}_{F}=\dfrac{f_{n}f_{n-1}\cdots f_{n-k+1}}{f_{k}f_{k-1}\cdots
f_{1}}.
\]
This is called the $\left(  n,k\right)  $-th \textit{Fibonomial coefficient}
(in analogy to the binomial coefficient $\dbinom{n}{k}=\dfrac{n\left(
n-1\right)  \cdots\left(  n-k+1\right)  }{k\left(  k-1\right)  \cdots1}$).

Show that $\dbinom{n}{k}_{F}$ is an integer.

\textbf{(b)} Try to extend as many identities for binomial coefficients as you
can to Fibonomial coefficients.

\textbf{(c)} Generalize to $\left(  a,b\right)  $-recurrent sequences with
arbitrary $a$ and $b$.
\end{exercise}

\section{\label{chp.perm}Permutations}

This chapter is devoted to permutations. We first recall how they are defined.

\subsection{Permutations and the symmetric group}

\begin{definition}
\label{def.composition}First, let us stipulate, once and for all, how we
define the composition of two maps: If $X$, $Y$ and $Z$ are three sets, and if
$\alpha:X\rightarrow Y$ and $\beta:Y\rightarrow Z$ are two maps, then
$\beta\circ\alpha$ denotes the map from $X$ to $Z$ which sends every $x\in X$
to $\beta\left(  \alpha\left(  x\right)  \right)  $. This map $\beta
\circ\alpha$ is called the \textit{composition} of $\beta$ and $\alpha$ (and
is sometimes abbreviated as $\beta\alpha$). This is the classical notation for
composition of maps, and the reason why I am so explicitly reminding you of it
is that some people (e.g., Herstein in \cite{Herstein}) use a different
convention that conflicts with it: They write maps \textquotedblleft on the
right\textquotedblright\ (i.e., they denote the image of an element $x\in X$
under the map $\alpha:X\rightarrow Y$ by $x^{\alpha}$ or $x\alpha$ instead of
$\alpha\left(  x\right)  $), and they define composition \textquotedblleft the
other way round\textquotedblright\ (i.e., they write $\alpha\circ\beta$ for
what we call $\beta\circ\alpha$). They have reasons for what they are doing,
but I shall use the classical notation because most of the literature agrees
with it.
\end{definition}

\begin{definition}
Let us also recall what it means for two maps to be \textit{inverse}.

Let $X$ and $Y$ be two sets. Two maps $f : X \to Y$ and $g : Y \to X$ are said
to be \textit{mutually inverse} if they satisfy $g \circ f = \operatorname{id}%
_{X}$ and $f \circ g = \operatorname{id}_{Y}$. (In other words, two maps $f :
X \to Y$ and $g : Y \to X$ are mutually inverse if and only if every $x \in X$
satisfies $g\left(  f\left(  x\right)  \right)  = x$ and every $y \in Y$
satisfies $f\left(  g\left(  y\right)  \right)  = y$.)

Let $f : X \to Y$ be a map. If there exists a map $g : Y \to X$ such that $f$
and $g$ are mutually inverse, then this map $g$ is unique (this is easy to
check) and is called the \textit{inverse} of $f$ and denoted by $f^{-1}$. In
this case, the map $f$ is said to be \textit{invertible}. It is easy to see
that if $g$ is the inverse of $f$, then $f$ is the inverse of $g$.

It is well-known that a map $f : X \to Y$ is invertible if and only if $f$ is
bijective (i.e., both injective and surjective). The words ``invertible'' and
``bijective'' are thus synonyms (at least when used for a map between two sets
-- in other situations, they can be rather different). Nevertheless, both of
them are commonly used, often by the same authors (since they convey slightly
different mental images).

A bijective map is also called a \textit{bijection} or a \textit{1-to-1
correspondence} (or a \textit{one-to-one correspondence}). When there is a
bijection from $X$ to $Y$, one says that the elements of $X$ are \textit{in
bijection with} (or \textit{in one-to-one correspondence with}) the elements
of $Y$. It is well-known that two sets $X$ and $Y$ have the same cardinality
if and only if there exists a bijection from $X$ to $Y$. (This is precisely
Theorem \ref{thm.bijection=eqsize}.)
\end{definition}

\begin{definition}
\label{def.permutation}A \textit{permutation} of a set $X$ means a bijection
from $X$ to $X$. The permutations of a given set $X$ can be composed (i.e., if
$\alpha$ and $\beta$ are two permutations of $X$, then so is $\alpha\circ
\beta$) and have inverses (which, again, are permutations of $X$). More precisely:

\begin{itemize}
\item If $\alpha$ and $\beta$ are two permutations of a given set $X$, then
the composition $\alpha\circ\beta$ is again a permutation of $X$.

\item Any three permutations $\alpha$, $\beta$ and $\gamma$ of $X$ satisfy
$\left(  \alpha\circ\beta\right)  \circ\gamma=\alpha\circ\left(  \beta
\circ\gamma\right)  $. (This holds, more generally, for arbitrary maps which
can be composed.)

\item The identity map $\operatorname*{id}:X\rightarrow X$ (this is the map
which sends every element $x\in X$ to itself) is a permutation of $X$; it is
also called the \textit{identity permutation}. Every permutation $\alpha$ of
$X$ satisfies $\operatorname*{id}\circ\alpha=\alpha$ and $\alpha
\circ\operatorname*{id}=\alpha$. (Again, this can be generalized to arbitrary maps.)

\item For every permutation $\alpha$ of $X$, the inverse map $\alpha^{-1}$ is
well-defined and is again a permutation of $X$. We have $\alpha\circ
\alpha^{-1}=\operatorname*{id}$ and $\alpha^{-1}\circ\alpha=\operatorname*{id}%
$.
\end{itemize}

In the lingo of algebraists, these four properties show that the set of all
permutations of $X$ is a
\href{https://en.wikipedia.org/?title=Group (mathematics)}{group} whose binary
operation is composition, and whose neutral element is the identity
permutation $\operatorname*{id}:X\rightarrow X$. This group is known as the
\textit{symmetric group of the set }$X$. (We will define the notion of a group
later, in Definition \ref{def.group}; thus you might not understand the
preceding two sentences at this point. If you do not care about groups, you
should just remember that the symmetric group of $X$ is the set of all
permutations of $X$.)
\end{definition}

\begin{remark}
Some authors define a permutation of a finite set $X$ to mean a list of all
elements of $X$, each occurring exactly once. This is \textbf{not} the meaning
that the word \textquotedblleft permutation\textquotedblright\ has in these
notes! It is a different notion which, for historical reasons, has been called
\textquotedblleft permutation\textquotedblright\ as well. On
\href{https://en.wikipedia.org/wiki/Permutation\%23Definition_and_one-line_notation}{the
Wikipedia page for \textquotedblleft permutation\textquotedblright}, the two
notions are called \textquotedblleft active\textquotedblright\ and
\textquotedblleft passive\textquotedblright, respectively: An
\textquotedblleft active\textquotedblright\ permutation of $X$ means a
bijection from $X$ to $X$ (that is, a permutation of $X$ in our meaning of
this word), whereas a \textquotedblleft passive\textquotedblright\ permutation
of $X$ means a list of all elements of $X$, each occurring exactly once. For
example, if $X=\left\{  \text{\textquotedblleft cat\textquotedblright,
\textquotedblleft dog\textquotedblright, \textquotedblleft
archaeopteryx\textquotedblright}\right\}  $, then the map%
\begin{align*}
\text{\textquotedblleft cat\textquotedblright\ }  &  \mapsto\text{
\textquotedblleft archaeopteryx\textquotedblright},\\
\text{\textquotedblleft archaeopteryx\textquotedblright\ }  &  \mapsto\text{
\textquotedblleft dog\textquotedblright},\\
\text{\textquotedblleft dog\textquotedblright\ }  &  \mapsto\text{
\textquotedblleft cat\textquotedblright}%
\end{align*}
is an \textquotedblleft active\textquotedblright\ permutation of $X$, whereas
the list $\left(  \text{\textquotedblleft dog\textquotedblright,
\textquotedblleft cat\textquotedblright, \textquotedblleft
archaeopteryx\textquotedblright}\right)  $ is a \textquotedblleft
passive\textquotedblright\ permutation of $X$.

When $X$ is the set $\left\{  1,2,\ldots,n\right\}  $ for some $n\in
\mathbb{N}$, then it is possible to equate each \textquotedblleft
active\textquotedblright\ permutation of $X$ with a \textquotedblleft
passive\textquotedblright\ permutation of $X$ (namely, its one-line notation,
defined below). More generally, this can be done when $X$ comes with a fixed
total order. In general, if $X$ is a finite set, then the number of
\textquotedblleft active\textquotedblright\ permutations of $X$ equals the
number of \textquotedblleft passive\textquotedblright\ permutations of $X$
(and both numbers equal $\left\vert X\right\vert !$), but until you fix some
ordering of the elements of $X$, there is no \textquotedblleft
natural\textquotedblright\ way to match the \textquotedblleft
passive\textquotedblright\ permutations with the \textquotedblleft
active\textquotedblright\ ones. (And when $X$ is infinite, the notion of a
\textquotedblleft passive\textquotedblright\ permutation is not even well-defined.)

To reiterate: For us, the word \textquotedblleft permutation\textquotedblright%
\ shall always mean an \textquotedblleft active\textquotedblright\ permutation!
\end{remark}

Recall that $\mathbb{N}=\left\{  0,1,2,\ldots\right\}  $.

\begin{definition}
Let $n\in\mathbb{N}$.

Let $S_{n}$ be the symmetric group of the set $\left\{  1,2,\ldots,n\right\}
$. This is the set of all permutations of the set $\left\{  1,2,\ldots
,n\right\}  $. It contains the identity permutation $\operatorname*{id}\in
S_{n}$ which sends every $i\in\left\{  1,2,\ldots,n\right\}  $ to $i$.
\end{definition}

A well-known fact states that for every $n\in\mathbb{N}$, the size of the
symmetric group $S_{n}$ is $\left\vert S_{n}\right\vert =n!$ (that is, there
are exactly $n!$ permutations of $\left\{  1,2,\ldots,n\right\}  $). (One
proof of this fact -- not the simplest -- is given in the proof of Corollary
\ref{cor.transpos.code.n!} below.)

We will often write a permutation $\sigma\in S_{n}$ as the list $\left(
\sigma\left(  1\right)  ,\sigma\left(  2\right)  ,\ldots,\sigma\left(
n\right)  \right)  $ of its values. This is known as the \textit{one-line
notation} for permutations (because it is a single-rowed list, as opposed to
e.g. the two-line notation which is a two-rowed
table).\footnote{Combinatorialists often omit the parentheses and the commas
(i.e., they just write $\sigma\left(  1\right)  \sigma\left(  2\right)
\cdots\sigma\left(  n\right)  $, hoping that noone will mistake this for a
product), since there is unfortunately another notation for permutations (the
\textit{cycle notation}) which also writes them as lists (actually, lists of
lists) but where the lists have a different meaning.} For instance, the
permutation in $S_{3}$ which sends $1$ to $2$, $2$ to $1$ and $3$ to $3$ is
written $\left(  2,1,3\right)  $ in one-line notation.

The exact relation between lists and permutations is given by the following
simple fact:

\begin{proposition}
\label{prop.perms.lists}Let $n\in\mathbb{N}$. Let $\left[  n\right]  =\left\{
1,2,\ldots,n\right\}  $.

\textbf{(a)} If $\sigma\in S_{n}$, then each element of $\left[  n\right]  $
appears exactly once in the list $\left(  \sigma\left(  1\right)
,\sigma\left(  2\right)  ,\ldots,\sigma\left(  n\right)  \right)  $.

\textbf{(b)} If $\left(  p_{1},p_{2},\ldots,p_{n}\right)  $ is a list of
elements of $\left[  n\right]  $ such that each element of $\left[  n\right]
$ appears exactly once in this list $\left(  p_{1},p_{2},\ldots,p_{n}\right)
$, then there exists a unique permutation $\sigma\in S_{n}$ such that $\left(
p_{1},p_{2},\ldots,p_{n}\right)  =\left(  \sigma\left(  1\right)
,\sigma\left(  2\right)  ,\ldots,\sigma\left(  n\right)  \right)  $.

\textbf{(c)} Let $k\in\left\{  0,1,\ldots,n\right\}  $. If $\left(
p_{1},p_{2},\ldots,p_{k}\right)  $ is a list of some elements of $\left[
n\right]  $ such that $p_{1},p_{2},\ldots,p_{k}$ are distinct, then there
exists a permutation $\sigma\in S_{n}$ such that $\left(  p_{1},p_{2}%
,\ldots,p_{k}\right)  =\left(  \sigma\left(  1\right)  ,\sigma\left(
2\right)  ,\ldots,\sigma\left(  k\right)  \right)  $.
\end{proposition}

At this point, let us clarify what we mean by \textquotedblleft
distinct\textquotedblright: Several objects $u_{1},u_{2},\ldots,u_{k}$ are
said to be \textit{distinct} if every $i\in\left\{  1,2,\ldots,k\right\}  $
and $j\in\left\{  1,2,\ldots,k\right\}  $ satisfying $i\neq j$ satisfy
$u_{i}\neq u_{j}$. (Some people call this \textquotedblleft pairwise
distinct\textquotedblright.) So, for example, the numbers $2,1,6$ are
distinct, but the numbers $6,1,6$ are not (although $6$ and $1$ are distinct).
Instead of saying that some objects $u_{1},u_{2},\ldots,u_{k}$ are distinct,
we can also say that \textquotedblleft the list $\left(  u_{1},u_{2}%
,\ldots,u_{k}\right)  $ has no repetitions\textquotedblright\footnote{A
repetition just means an element which occurs more than once in the list. It
does not matter whether the occurrences are at consecutive positions or not.}.

\begin{remark}
The $\sigma$ in Proposition \ref{prop.perms.lists} \textbf{(b)} is uniquely
determined, but the $\sigma$ in Proposition \ref{prop.perms.lists}
\textbf{(c)} is not (in general). More precisely, in Proposition
\ref{prop.perms.lists} \textbf{(c)}, there are $\left(  n-k\right)  !$
possible choices of $\sigma$ that work. (This is easy to check.)
\end{remark}

\begin{proof}
[Proof of Proposition \ref{prop.perms.lists}.]Proposition
\ref{prop.perms.lists} is a basic fact and its proof is simple. I am going to
present the proof in great detail, but you are not missing much if you skip it
for its obviousness (just make sure you know \textbf{why} it is obvious).

Recall that $S_{n}$ is the set of all permutations of the set $\left\{
1,2,\ldots,n\right\}  $. In other words, $S_{n}$ is the set of all
permutations of the set $\left[  n\right]  $ (since $\left\{  1,2,\ldots
,n\right\}  =\left[  n\right]  $).

\textbf{(a)} Let $\sigma\in S_{n}$. Let $i\in\left[  n\right]  $.

We have $\sigma\in S_{n}$. In other words, $\sigma$ is a permutation of
$\left[  n\right]  $ (since $S_{n}$ is the set of all permutations of the set
$\left[  n\right]  $). In other words, $\sigma$ is a bijective map $\left[
n\right]  \rightarrow\left[  n\right]  $. Hence, $\sigma$ is both surjective
and injective.

Now, we make the following two observations:

\begin{itemize}
\item The number $i$ appears in the list $\left(  \sigma\left(  1\right)
,\sigma\left(  2\right)  ,\ldots,\sigma\left(  n\right)  \right)
$\ \ \ \ \footnote{\textit{Proof.} The map $\sigma$ is surjective. Hence,
there exists some $j\in\left[  n\right]  $ such that $i=\sigma\left(
j\right)  $. In other words, the number $i$ appears in the list $\left(
\sigma\left(  1\right)  ,\sigma\left(  2\right)  ,\ldots,\sigma\left(
n\right)  \right)  $. Qed.}.

\item The number $i$ appears at most once in the list $\left(  \sigma\left(
1\right)  ,\sigma\left(  2\right)  ,\ldots,\sigma\left(  n\right)  \right)
$\ \ \ \ \footnote{\textit{Proof.} Let us assume the contrary (for the sake of
contradiction). Thus, $i$ appears more than once in the list $\left(
\sigma\left(  1\right)  ,\sigma\left(  2\right)  ,\ldots,\sigma\left(
n\right)  \right)  $. In other words, $i$ appears at least twice in this list.
In other words, there exist two distinct elements $p$ and $q$ of $\left[
n\right]  $ such that $\sigma\left(  p\right)  =i$ and $\sigma\left(
q\right)  =i$. Consider these $p$ and $q$.
\par
We have $p\neq q$ (since $p$ and $q$ are distinct), so that $\sigma\left(
p\right)  \neq\sigma\left(  q\right)  $ (since $\sigma$ is injective). This
contradicts $\sigma\left(  p\right)  =i=\sigma\left(  q\right)  $. This
contradiction proves that our assumption was wrong, qed.}.
\end{itemize}

Combining these two observations, we conclude that the number $i$ appears
exactly once in the list $\left(  \sigma\left(  1\right)  ,\sigma\left(
2\right)  ,\ldots,\sigma\left(  n\right)  \right)  $.

Let us now forget that we fixed $i$. We thus have shown that if $i\in\left[
n\right]  $, then $i$ appears exactly once in the list $\left(  \sigma\left(
1\right)  ,\sigma\left(  2\right)  ,\ldots,\sigma\left(  n\right)  \right)  $.
In other words, each element of $\left[  n\right]  $ appears exactly once in
the list $\left(  \sigma\left(  1\right)  ,\sigma\left(  2\right)
,\ldots,\sigma\left(  n\right)  \right)  $. This proves Proposition
\ref{prop.perms.lists} \textbf{(a)}.

\textbf{(b)} Let $\left(  p_{1},p_{2},\ldots,p_{n}\right)  $ be a list of
elements of $\left[  n\right]  $ such that each element of $\left[  n\right]
$ appears exactly once in this list $\left(  p_{1},p_{2},\ldots,p_{n}\right)
$.

We have $p_{i}\in\left[  n\right]  $ for every $i\in\left[  n\right]  $ (since
$\left(  p_{1},p_{2},\ldots,p_{n}\right)  $ is a list of elements of $\left[
n\right]  $).

We define a map $\tau:\left[  n\right]  \rightarrow\left[  n\right]  $ by
setting%
\begin{equation}
\left(  \tau\left(  i\right)  =p_{i}\ \ \ \ \ \ \ \ \ \ \text{for every }%
i\in\left[  n\right]  \right)  . \label{pf.prop.perms.lists.b.1}%
\end{equation}
(This is well-defined, because we have $p_{i}\in\left[  n\right]  $ for every
$i\in\left[  n\right]  $.) The map $\tau$ is
injective\footnote{\textit{Proof.} Let $u$ and $v$ be two elements of $\left[
n\right]  $ such that $\tau\left(  u\right)  =\tau\left(  v\right)  $. We
shall show that $u=v$.
\par
Indeed, we assume the contrary (for the sake of contradiction). Thus, $u\neq
v$.
\par
The definition of $\tau\left(  u\right)  $ shows that $\tau\left(  u\right)
=p_{u}$. But we also have $\tau\left(  u\right)  =\tau\left(  v\right)
=p_{v}$ (by the definition of $\tau\left(  v\right)  $). Now, the element
$\tau\left(  u\right)  $ of $\left[  n\right]  $ appears (at least) twice in
the list $\left(  p_{1},p_{2},\ldots,p_{n}\right)  $: once at the $u$-th
position (since $\tau\left(  u\right)  =p_{u}$), and again at the $v$-th
position (since $\tau\left(  u\right)  =p_{v}$). (And these are two distinct
positions, because $u\neq v$.)
\par
But let us recall that each element of $\left[  n\right]  $ appears exactly
once in this list $\left(  p_{1},p_{2},\ldots,p_{n}\right)  $. Hence, no
element of $\left[  n\right]  $ appears more than once in the list $\left(
p_{1},p_{2},\ldots,p_{n}\right)  $. In particular, $\tau\left(  u\right)  $
cannot appear more than once in this list $\left(  p_{1},p_{2},\ldots
,p_{n}\right)  $. This contradicts the fact that $\tau\left(  u\right)  $
appears twice in the list $\left(  p_{1},p_{2},\ldots,p_{n}\right)  $.
\par
This contradiction shows that our assumption was wrong. Hence, $u=v$ is
proven.
\par
Now, let us forget that we fixed $u$ and $v$. We thus have proven that if $u$
and $v$ are two elements of $\left[  n\right]  $ such that $\tau\left(
u\right)  =\tau\left(  v\right)  $, then $u=v$. In other words, the map $\tau$
is injective. Qed.} and surjective\footnote{\textit{Proof.} Let $u\in\left[
n\right]  $. Each element of $\left[  n\right]  $ appears exactly once in the
list $\left(  p_{1},p_{2},\ldots,p_{n}\right)  $. Applying this to the element
$u$ of $\left[  n\right]  $, we conclude that $u$ appears exactly once in the
list $\left(  p_{1},p_{2},\ldots,p_{n}\right)  $. In other words, there exists
exactly one $i\in\left[  n\right]  $ such that $u=p_{i}$. Consider this $i$.
The definition of $\tau$ yields $\tau\left(  i\right)  =p_{i}$. Compared with
$u=p_{i}$, this yields $\tau\left(  i\right)  =u$.
\par
Hence, there exists a $j\in\left[  n\right]  $ such that $\tau\left(
j\right)  =u$ (namely, $j=i$).
\par
Let us now forget that we fixed $u$. We thus have proven that for every
$u\in\left[  n\right]  $, there exists a $j\in\left[  n\right]  $ such that
$\tau\left(  j\right)  =u$. In other words, the map $\tau$ is surjective.
Qed.}. Hence, the map $\tau$ is bijective. In other words, $\tau$ is a
permutation of $\left[  n\right]  $ (since $\tau$ is a map $\left[  n\right]
\rightarrow\left[  n\right]  $). In other words, $\tau\in S_{n}$ (since
$S_{n}$ is the set of all permutations of the set $\left[  n\right]  $).
Clearly, $\left(  \tau\left(  1\right)  ,\tau\left(  2\right)  ,\ldots
,\tau\left(  n\right)  \right)  =\left(  p_{1},p_{2},\ldots,p_{n}\right)  $
(because of (\ref{pf.prop.perms.lists.b.1})), so that $\left(  p_{1}%
,p_{2},\ldots,p_{n}\right)  =\left(  \tau\left(  1\right)  ,\tau\left(
2\right)  ,\ldots,\tau\left(  n\right)  \right)  $.

Hence, there exists a permutation $\sigma\in S_{n}$ such that \newline$\left(
p_{1},p_{2},\ldots,p_{n}\right)  =\left(  \sigma\left(  1\right)
,\sigma\left(  2\right)  ,\ldots,\sigma\left(  n\right)  \right)  $ (namely,
$\sigma=\tau$). Moreover, there exists \textbf{at most one} such
permutation\footnote{\textit{Proof.} Let $\sigma_{1}$ and $\sigma_{2}$ be two
permutations $\sigma\in S_{n}$ such that $\left(  p_{1},p_{2},\ldots
,p_{n}\right)  =\left(  \sigma\left(  1\right)  ,\sigma\left(  2\right)
,\ldots,\sigma\left(  n\right)  \right)  $. Thus, $\sigma_{1}$ is a
permutation in $S_{n}$ such that $\left(  p_{1},p_{2},\ldots,p_{n}\right)
=\left(  \sigma_{1}\left(  1\right)  ,\sigma_{1}\left(  2\right)
,\ldots,\sigma_{1}\left(  n\right)  \right)  $, and $\sigma_{2}$ is a
permutation in $S_{n}$ such that $\left(  p_{1},p_{2},\ldots,p_{n}\right)
=\left(  \sigma_{2}\left(  1\right)  ,\sigma_{2}\left(  2\right)
,\ldots,\sigma_{2}\left(  n\right)  \right)  $.
\par
We have $\left(  \sigma_{1}\left(  1\right)  ,\sigma_{1}\left(  2\right)
,\ldots,\sigma_{1}\left(  n\right)  \right)  =\left(  p_{1},p_{2},\ldots
,p_{n}\right)  =\left(  \sigma_{2}\left(  1\right)  ,\sigma_{2}\left(
2\right)  ,\ldots,\sigma_{2}\left(  n\right)  \right)  $. In other words,
every $i\in\left[  n\right]  $ satisfies $\sigma_{1}\left(  i\right)
=\sigma_{2}\left(  i\right)  $. In other words, $\sigma_{1}=\sigma_{2}$.
\par
Let us now forget that we fixed $\sigma_{1}$ and $\sigma_{2}$. We thus have
shown that if $\sigma_{1}$ and $\sigma_{2}$ are two permutations $\sigma\in
S_{n}$ such that $\left(  p_{1},p_{2},\ldots,p_{n}\right)  =\left(
\sigma\left(  1\right)  ,\sigma\left(  2\right)  ,\ldots,\sigma\left(
n\right)  \right)  $, then $\sigma_{1}=\sigma_{2}$. In other words, any two
permutations $\sigma\in S_{n}$ such that $\left(  p_{1},p_{2},\ldots
,p_{n}\right)  =\left(  \sigma\left(  1\right)  ,\sigma\left(  2\right)
,\ldots,\sigma\left(  n\right)  \right)  $ must be equal to each other. In
other words, there exists \textbf{at most one} permutation $\sigma\in S_{n}$
such that $\left(  p_{1},p_{2},\ldots,p_{n}\right)  =\left(  \sigma\left(
1\right)  ,\sigma\left(  2\right)  ,\ldots,\sigma\left(  n\right)  \right)  $.
Qed.}. Combining the claims of the previous two sentences, we conclude that
there exists a unique permutation $\sigma\in S_{n}$ such that $\left(
p_{1},p_{2},\ldots,p_{n}\right)  =\left(  \sigma\left(  1\right)
,\sigma\left(  2\right)  ,\ldots,\sigma\left(  n\right)  \right)  $. This
proves Proposition \ref{prop.perms.lists} \textbf{(b)}.

\textbf{(c)} Let $\left(  p_{1},p_{2},\ldots,p_{k}\right)  $ be a list of some
elements of $\left[  n\right]  $ such that $p_{1},p_{2},\ldots,p_{k}$ are
distinct. Thus, the list $\left(  p_{1},p_{2},\ldots,p_{k}\right)  $ contains
$k$ of the $n$ elements of $\left[  n\right]  $ (because $p_{1},p_{2}%
,\ldots,p_{k}$ are distinct). Let $q_{1},q_{2},\ldots,q_{n-k}$ be the
remaining $n-k$ elements of $\left[  n\right]  $ (listed in any arbitrary
order, with no repetition). Then, $\left(  p_{1},p_{2},\ldots,p_{k}%
,q_{1},q_{2},\ldots,q_{n-k}\right)  $ is a list of all $n$ elements of
$\left[  n\right]  $, with no repetitions\footnote{It has no repetitions
because:
\par
\begin{itemize}
\item there are no repetitions among $p_{1},p_{2},\ldots,p_{k}$;
\par
\item there are no repetitions among $q_{1},q_{2},\ldots,q_{n-k}$;
\par
\item the two lists $\left(  p_{1},p_{2},\ldots,p_{k}\right)  $ and $\left(
q_{1},q_{2},\ldots,q_{n-k}\right)  $ have no elements in common (because we
defined $q_{1},q_{2},\ldots,q_{n-k}$ to be the \textquotedblleft
remaining\textquotedblright\ $n-k$ elements of $\left[  n\right]  $, where
\textquotedblleft remaining\textquotedblright\ means \textquotedblleft not
contained in the list $\left(  p_{1},p_{2},\ldots,p_{k}\right)  $%
\textquotedblright).
\end{itemize}
}. In other words, each element of $\left[  n\right]  $ appears exactly once
in this list $\left(  p_{1},p_{2},\ldots,p_{k},q_{1},q_{2},\ldots
,q_{n-k}\right)  $ (and each entry in this list is an element of $\left[
n\right]  $). Hence, we can apply Proposition \ref{prop.perms.lists}
\textbf{(b)} to $\left(  p_{1},p_{2},\ldots,p_{k},q_{1},q_{2},\ldots
,q_{n-k}\right)  $ instead of $\left(  p_{1},p_{2},\ldots,p_{n}\right)  $. As
a consequence, we conclude that there exists a unique permutation $\sigma\in
S_{n}$ such that $\left(  p_{1},p_{2},\ldots,p_{k},q_{1},q_{2},\ldots
,q_{n-k}\right)  =\left(  \sigma\left(  1\right)  ,\sigma\left(  2\right)
,\ldots,\sigma\left(  n\right)  \right)  $. Let $\tau$ be this $\sigma$.

Thus, $\tau\in S_{n}$ is a permutation such that
\[
\left(  p_{1},p_{2},\ldots,p_{k},q_{1},q_{2},\ldots,q_{n-k}\right)  =\left(
\tau\left(  1\right)  ,\tau\left(  2\right)  ,\ldots,\tau\left(  n\right)
\right)  .
\]
Now,%
\begin{align*}
&  \left(  p_{1},p_{2},\ldots,p_{k}\right) \\
&  =\left(  \text{the list of the first }k\text{ entries of the list
}\underbrace{\left(  p_{1},p_{2},\ldots,p_{k},q_{1},q_{2},\ldots
,q_{n-k}\right)  }_{=\left(  \tau\left(  1\right)  ,\tau\left(  2\right)
,\ldots,\tau\left(  n\right)  \right)  }\right) \\
&  =\left(  \text{the list of the first }k\text{ entries of the list }\left(
\tau\left(  1\right)  ,\tau\left(  2\right)  ,\ldots,\tau\left(  n\right)
\right)  \right) \\
&  =\left(  \tau\left(  1\right)  ,\tau\left(  2\right)  ,\ldots,\tau\left(
k\right)  \right)  .
\end{align*}
Hence, there exists a permutation $\sigma\in S_{n}$ such that \newline$\left(
p_{1},p_{2},\ldots,p_{k}\right)  =\left(  \sigma\left(  1\right)
,\sigma\left(  2\right)  ,\ldots,\sigma\left(  k\right)  \right)  $ (namely,
$\sigma=\tau$). This proves Proposition \ref{prop.perms.lists} \textbf{(c)}.
\end{proof}

\subsection{Inversions, lengths and the permutations $s_{i} \in S_{n}$}

\begin{definition}
\label{def.perm.si}Let $n\in\mathbb{N}$. For each $i\in\left\{  1,2,\ldots
,n-1\right\}  $, let $s_{i}$ be the permutation in $S_{n}$ that swaps $i$ with
$i+1$ but leaves all other numbers unchanged. Formally speaking, $s_{i}$ is
the permutation in $S_{n}$ given by%
\[
\left(  s_{i}\left(  k\right)  =%
\begin{cases}
i+1, & \text{if }k=i;\\
i, & \text{if }k=i+1;\\
k, & \text{if }k\notin\left\{  i,i+1\right\}
\end{cases}
\ \ \ \ \ \ \ \ \ \ \text{for all }k\in\left\{  1,2,\ldots,n\right\}  \right)
.
\]
Thus, in one-line notation%
\[
s_{i}=\left(  1,2,\ldots,i-1,i+1,i,i+2,\ldots,n\right)  .
\]
Notice that $s_{i}^{2}=\operatorname*{id}$ for every $i\in\left\{
1,2,\ldots,n-1\right\}  $. (Here, we are using the notation $\alpha^{2}$ for
$\alpha\circ\alpha$, where $\alpha$ is a permutation in $S_{n}$.)
\end{definition}

\begin{verlong}
\begin{proof}
[Proof of $s_{i}^{2}=\operatorname*{id}$ for every $i\in\left\{
1,2,\ldots,n-1\right\}  $.]For the sake of completeness, let us prove the
preceding claim (that is, $s_{i}^{2}=\operatorname*{id}$ for every
$i\in\left\{  1,2,\ldots,n-1\right\}  $) as well as the fact that $s_{i}$ is a
permutation whenever $i\in\left\{  1,2,\ldots,n-1\right\}  $. We shall prove
both of these facts in one swoop:

\begin{statement}
\textit{Claim 1:} Let $i\in\left\{  1,2,\ldots,n\right\}  $. Then, $s_{i}$ is
a permutation in $S_{n}$ and satisfies $s_{i}^{2}=\operatorname*{id}$.
\end{statement}

[\textit{Proof of Claim 1:} Let $i\in\left\{  1,2,\ldots,n-1\right\}  $. Thus,
$i+1\in\left\{  2,3,\ldots,n\right\}  \subseteq\left\{  1,2,\ldots,n\right\}
$ and $i\in\left\{  1,2,\ldots,n-1\right\}  \subseteq\left\{  1,2,\ldots
,n\right\}  $. Now, for each $k\in\left\{  1,2,\ldots,n\right\}  $, we have%
\[%
\begin{cases}
i+1, & \text{if }k=i;\\
i, & \text{if }k=i+1;\\
k, & \text{if }k\notin\left\{  i,i+1\right\}
\end{cases}
\in\left\{  1,2,\ldots,n\right\}
\]
(since all three possible values $i+1$, $i$ and $k$ of the expression $%
\begin{cases}
i+1, & \text{if }k=i;\\
i, & \text{if }k=i+1;\\
k, & \text{if }k\notin\left\{  i,i+1\right\}
\end{cases}
$ belong to $\left\{  1,2,\ldots,n\right\}  $). Thus, the map $s_{i}:\left\{
1,2,\ldots,n\right\}  \rightarrow\left\{  1,2,\ldots,n\right\}  $ is well-defined.

Next, we will prove that $s_{i}^{2}=\operatorname*{id}$.

Fix $k\in\left\{  1,2,\ldots,n\right\}  $. We shall show that $s_{i}%
^{2}\left(  k\right)  =k$.

Indeed, we are in one of the following three cases:

\textit{Case 1:} We have $k=i$.

\textit{Case 2:} We have $k=i+1$.

\textit{Case 3:} We have neither $k=i$ nor $k=i+1$.

Let us first consider Case 1. In this case, we have $k=i$. But the definition
of $s_{i}$ yields%
\[
s_{i}\left(  k\right)  =%
\begin{cases}
i+1, & \text{if }k=i;\\
i, & \text{if }k=i+1;\\
k, & \text{if }k\notin\left\{  i,i+1\right\}
\end{cases}
=i+1\ \ \ \ \ \ \ \ \ \ \left(  \text{since }k=i\right)  .
\]
Now,%
\begin{align*}
\underbrace{s_{i}^{2}}_{=s_{i}\circ s_{i}}\left(  k\right)   &  =\left(
s_{i}\circ s_{i}\right)  \left(  k\right)  =s_{i}\left(  \underbrace{s_{i}%
\left(  k\right)  }_{=i+1}\right)  =s_{i}\left(  i+1\right) \\
&  =%
\begin{cases}
i+1, & \text{if }i+1=i;\\
i, & \text{if }i+1=i+1;\\
i+1, & \text{if }i+1\notin\left\{  i,i+1\right\}
\end{cases}
\ \ \ \ \ \ \ \ \ \ \left(  \text{by the definition of }s_{i}\right) \\
&  =i\ \ \ \ \ \ \ \ \ \ \left(  \text{since }i+1=i+1\right) \\
&  =k\ \ \ \ \ \ \ \ \ \ \left(  \text{since }k=i\right)  .
\end{align*}
Hence, $s_{i}^{2}\left(  k\right)  =k$ is proven in Case 1.

Let us next consider Case 2. In this case, we have $k=i+1$. But the definition
of $s_{i}$ yields%
\[
s_{i}\left(  k\right)  =%
\begin{cases}
i+1, & \text{if }k=i;\\
i, & \text{if }k=i+1;\\
k, & \text{if }k\notin\left\{  i,i+1\right\}
\end{cases}
=i\ \ \ \ \ \ \ \ \ \ \left(  \text{since }k=i+1\right)  .
\]
Now,%
\begin{align*}
\underbrace{s_{i}^{2}}_{=s_{i}\circ s_{i}}\left(  k\right)   &  =\left(
s_{i}\circ s_{i}\right)  \left(  k\right)  =s_{i}\left(  \underbrace{s_{i}%
\left(  k\right)  }_{=i}\right)  =s_{i}\left(  i\right) \\
&  =%
\begin{cases}
i+1, & \text{if }i=i;\\
i, & \text{if }i=i+1;\\
i, & \text{if }i\notin\left\{  i,i+1\right\}
\end{cases}
\ \ \ \ \ \ \ \ \ \ \left(  \text{by the definition of }s_{i}\right) \\
&  =i+1\ \ \ \ \ \ \ \ \ \ \left(  \text{since }i=i\right) \\
&  =k\ \ \ \ \ \ \ \ \ \ \left(  \text{since }k=i+1\right)  .
\end{align*}
Hence, $s_{i}^{2}\left(  k\right)  =k$ is proven in Case 2.

Let us finally consider Case 3. In this case, we have neither $k=i$ nor
$k=i+1$. In other words, $k$ equals none of the two numbers $i$ and $i+1$. In
other words, $k\notin\left\{  i,i+1\right\}  $. But the definition of $s_{i}$
yields%
\[
s_{i}\left(  k\right)  =%
\begin{cases}
i+1, & \text{if }k=i;\\
i, & \text{if }k=i+1;\\
k, & \text{if }k\notin\left\{  i,i+1\right\}
\end{cases}
=k\ \ \ \ \ \ \ \ \ \ \left(  \text{since }k\notin\left\{  i,i+1\right\}
\right)  .
\]
Now,%
\[
\underbrace{s_{i}^{2}}_{=s_{i}\circ s_{i}}\left(  k\right)  =\left(
s_{i}\circ s_{i}\right)  \left(  k\right)  =s_{i}\left(  \underbrace{s_{i}%
\left(  k\right)  }_{=k}\right)  =s_{i}\left(  k\right)  =k.
\]
Hence, $s_{i}^{2}\left(  k\right)  =k$ is proven in Case 3.

We have now proven $s_{i}^{2}\left(  k\right)  =k$ in each of the three Cases
1, 2 and 3. Since these three Cases cover all possibilities, we thus conclude
that $s_{i}^{2}\left(  k\right)  =k$ always holds. Thus, $s_{i}^{2}\left(
k\right)  =k=\operatorname*{id}\left(  k\right)  $.

Now, forget that we fixed $k$. We thus have shown that $s_{i}^{2}\left(
k\right)  =\operatorname*{id}\left(  k\right)  $ for each $k\in\left\{
1,2,\ldots,n\right\}  $. In other words, $s_{i}^{2}=\operatorname*{id}$ (since
both $s_{i}^{2}$ and $\operatorname*{id}$ are maps from $\left\{
1,2,\ldots,n\right\}  $ to $\left\{  1,2,\ldots,n\right\}  $). Hence,
$s_{i}\circ s_{i}=s_{i}^{2}=\operatorname*{id}$; thus, the maps $s_{i}$ and
$s_{i}$ are mutually inverse (since $s_{i}\circ s_{i}=\operatorname*{id}$ and
$s_{i}\circ s_{i}=\operatorname*{id}$). Hence, the map $s_{i}$ is invertible,
i.e., bijective. Thus, $s_{i}$ is a bijection from $\left\{  1,2,\ldots
,n\right\}  $ to $\left\{  1,2,\ldots,n\right\}  $. In other words, $s_{i}$ is
a permutation of $\left\{  1,2,\ldots,n\right\}  $ (by the definition of a
permutation). In other words, $s_{i}\in S_{n}$.

Altogether, we thus have shown that $s_{i}$ is a permutation in $S_{n}$ and
satisfies $s_{i}^{2}=\operatorname*{id}$. This proves Claim 1.]
\end{proof}
\end{verlong}

\begin{exercise}
\label{exe.ps2.2.4}Let $n\in\mathbb{N}$.

\textbf{(a)} Show that $s_{i}\circ s_{i+1}\circ s_{i}=s_{i+1}\circ s_{i}\circ
s_{i+1}$ for all $i\in\left\{  1,2,\ldots,n-2\right\}  $.

\textbf{(b)} Show that every permutation $\sigma\in S_{n}$ can be written as a
composition of several permutations of the form $s_{k}$ (with $k\in\left\{
1,2,\ldots,n-1\right\}  $). For example, if $n=3$, then the
permutation\footnotemark\ $\left(  3,1,2\right)  $ in $S_{3}$ can be written
as the composition $s_{2}\circ s_{1}$, while the permutation $\left(
3,2,1\right)  $ in $S_{3}$ can be written as the composition $s_{1}\circ
s_{2}\circ s_{1}$ or also as the composition $s_{2}\circ s_{1}\circ s_{2}$.

[\textbf{Hint:} If you do not immediately see why this works, consider reading further.]

\textbf{(c)} Let $w_{0}$ denote the permutation in $S_{n}$ which sends each
$k\in\left\{  1,2,\ldots,n\right\}  $ to $n+1-k$. (In one-line notation, this
$w_{0}$ is written as $\left(  n,n-1,\ldots,1\right)  $.) Find an
\textbf{explicit} way to write $w_{0}$ as a composition of several
permutations of the form $s_{i}$ (with $i\in\left\{  1,2,\ldots,n-1\right\}  $).
\end{exercise}

\footnotetext{Recall that we are writing permutations in one-line notation.
Thus, \textquotedblleft the permutation $\left(  3,1,2\right)  $ in $S_{3}%
$\textquotedblright\ means the permutation $\sigma\in S_{3}$ satisfying
$\left(  \sigma\left(  1\right)  ,\sigma\left(  2\right)  ,\sigma\left(
3\right)  \right)  =\left(  3,1,2\right)  $.}

\begin{remark}
Symmetric groups appear in almost all parts of mathematics; unsurprisingly,
there is no universally accepted notation for them. We are using the notation
$S_{n}$ for the $n$-th symmetric group; other common notations for it are
$\mathfrak{S}_{n}$, $\Sigma_{n}$ and $\operatorname*{Sym}\left(  n\right)  $.
The permutations that we call $s_{1},s_{2},\ldots,s_{n-1}$ are often called
$\sigma_{1},\sigma_{2},\ldots,\sigma_{n-1}$. As already mentioned in
Definition \ref{def.composition}, some people write the composition of maps
\textquotedblleft backwards\textquotedblright, which causes their $\sigma
\circ\tau$ to be our $\tau\circ\sigma$, etc.. (Sadly, most authors are so sure
that their notation is standard that they never bother to define it.)

In the language of group theory, the statement of Exercise \ref{exe.ps2.2.4}
\textbf{(b)} says (or, more precisely, yields) that the permutations
$s_{1},s_{2},\ldots,s_{n-1}$ generate the group $S_{n}$.
\end{remark}

\begin{definition}
\label{def.perm.inversions}Let $n\in\mathbb{N}$. Let $\sigma\in S_{n}$ be a permutation.

\textbf{(a)} An \textit{inversion} of $\sigma$ means a pair $\left(
i,j\right)  $ of integers satisfying $1\leq i<j\leq n$ and $\sigma\left(
i\right)  >\sigma\left(  j\right)  $. For instance, the inversions of the
permutation $\left(  3,1,2\right)  $ (again, shown here in one-line notation)
in $S_{3}$ are $\left(  1,2\right)  $ and $\left(  1,3\right)  $ (because
$3>1$ and $3>2$), while the only inversion of the permutation $\left(
1,3,2\right)  $ in $S_{3}$ is $\left(  2,3\right)  $ (since $3>2$).

\textbf{(b)} The \textit{length} of $\sigma$ means the number of inversions of
$\sigma$. This length is denoted by $\ell\left(  \sigma\right)  $; it is a
nonnegative integer.
\end{definition}

If $n\in\mathbb{N}$, then any $\sigma\in S_{n}$ satisfies $0\leq\ell\left(
\sigma\right)  \leq\dbinom{n}{2}$ (since the number of inversions of $\sigma$
is clearly no larger than the total number of pairs $\left(  i,j\right)  $ of
integers satisfying $1\leq i<j\leq n$; but the latter number is $\dbinom{n}%
{2}$). The only permutation in $S_{n}$ having length $0$ is the identity
permutation $\operatorname*{id}=\left(  1,2,\ldots,n\right)  \in S_{n}%
$\ \ \ \ \footnote{The fact that the identity permutation $\operatorname*{id}%
\in S_{n}$ has length $\ell\left(  \operatorname*{id}\right)  =0$ is trivial.
The fact that it is the only one such permutation is easy (it essentially
follows from Exercise \ref{exe.ps2.2.5} \textbf{(d)}).}.

\begin{exercise}
\label{exe.ps2.2.5}Let $n\in\mathbb{N}$.

\textbf{(a)} Show that every permutation $\sigma\in S_{n}$ and every
$k\in\left\{  1,2,\ldots,n-1\right\}  $ satisfy%
\begin{equation}
\ell\left(  \sigma\circ s_{k}\right)  =%
\begin{cases}
\ell\left(  \sigma\right)  +1, & \text{if }\sigma\left(  k\right)
<\sigma\left(  k+1\right)  ;\\
\ell\left(  \sigma\right)  -1, & \text{if }\sigma\left(  k\right)
>\sigma\left(  k+1\right)
\end{cases}
\label{eq.exe.2.5.a.1}%
\end{equation}
and%
\begin{equation}
\ell\left(  s_{k}\circ\sigma\right)  =%
\begin{cases}
\ell\left(  \sigma\right)  +1, & \text{if }\sigma^{-1}\left(  k\right)
<\sigma^{-1}\left(  k+1\right)  ;\\
\ell\left(  \sigma\right)  -1, & \text{if }\sigma^{-1}\left(  k\right)
>\sigma^{-1}\left(  k+1\right)
\end{cases}
. \label{eq.exe.2.5.a.2}%
\end{equation}

\textbf{(b)} Show that any two permutations $\sigma$ and $\tau$ in $S_{n}$
satisfy $\ell\left(  \sigma\circ\tau\right)  \equiv\ell\left(  \sigma\right)
+\ell\left(  \tau\right)  \operatorname{mod}2$.

\textbf{(c)} Show that any two permutations $\sigma$ and $\tau$ in $S_{n}$
satisfy $\ell\left(  \sigma\circ\tau\right)  \leq\ell\left(  \sigma\right)
+\ell\left(  \tau\right)  $.

\textbf{(d)} If $\sigma\in S_{n}$ is a permutation satisfying $\sigma\left(
1\right)  \leq\sigma\left(  2\right)  \leq\cdots\leq\sigma\left(  n\right)  $,
then show that $\sigma=\operatorname*{id}$.

\textbf{(e)} Let $\sigma\in S_{n}$. Show that $\sigma$ can be written as a
composition of $\ell\left(  \sigma\right)  $ permutations of the form $s_{k}$
(with $k\in\left\{  1,2,\ldots,n-1\right\}  $).

\textbf{(f)} Let $\sigma\in S_{n}$. Then, show that $\ell\left(
\sigma\right)  =\ell\left(  \sigma^{-1}\right)  $.

\textbf{(g)} Let $\sigma\in S_{n}$. Show that $\ell\left(  \sigma\right)  $ is
the smallest $N\in\mathbb{N}$ such that $\sigma$ can be written as a
composition of $N$ permutations of the form $s_{k}$ (with $k\in\left\{
1,2,\ldots,n-1\right\}  $).
\end{exercise}

\begin{example}
\label{exa.2.5}Let us justify Exercise \ref{exe.ps2.2.5} \textbf{(a)} on an
example. The solution to Exercise \ref{exe.ps2.2.5} \textbf{(a)} given below
is essentially a (tiresome) formalization of the ideas seen in this example.

Let $n=5$, $k=3$ and $\sigma=\left(  4,2,1,5,3\right)  $ (written in one-line
notation). Then, $\sigma\circ s_{k}=\left(  4,2,5,1,3\right)  $; this is the
permutation obtained by swapping the $k$-th and the $\left(  k+1\right)  $-th
entry of $\sigma$ (where the word \textquotedblleft entry\textquotedblright%
\ refers to the one-line notation). On the other hand, $s_{k}\circ
\sigma=\left(  3,2,1,5,4\right)  $; this is the permutation obtained by
swapping the entry $k$ with the entry $k+1$ of $\sigma$. Mind the difference
between these two operations.

The inversions of $\sigma=\left(  4,2,1,5,3\right)  $ are $\left(  1,2\right)
$, $\left(  1,3\right)  $, $\left(  1,5\right)  $, $\left(  2,3\right)  $ and
$\left(  4,5\right)  $. These are the pairs $\left(  i,j\right)  $ of
positions such that $i$ is before $j$ (that is, $i<j$) but the $i$-th entry of
$\sigma$ is larger than the $j$-th entry of $\sigma$ (that is, $\sigma\left(
i\right)  >\sigma\left(  j\right)  $). In other words, these are the pairs of
positions at which the entries of $\sigma$ are out of order. On the other
hand, the inversions of $s_{k}\circ\sigma=\left(  3,2,1,5,4\right)  $ are
$\left(  1,2\right)  $, $\left(  1,3\right)  $, $\left(  2,3\right)  $ and
$\left(  4,5\right)  $. These are precisely the inversions of $\sigma$ except
for $\left(  1,5\right)  $. This is no surprise: In fact, $s_{k}\circ\sigma$
is obtained from $\sigma$ by swapping the entry $k$ with the entry $k+1$, and
this operation clearly preserves all inversions other than the one that is
directly being turned around (i.e., the inversion $\left(  i,j\right)  $ where
$\left\{  \sigma\left(  i\right)  ,\sigma\left(  j\right)  \right\}  =\left\{
k,k+1\right\}  $; in our case, this is the inversion $\left(  1,5\right)  $).
In general, when $\sigma^{-1}\left(  k\right)  >\sigma^{-1}\left(  k+1\right)
$ (that is, when $k$ appears further left than $k+1$ in the one-line notation
of $\sigma$), the inversions of $s_{k}\circ\sigma$ are the inversions of
$\sigma$ except for $\left(  \sigma^{-1}\left(  k+1\right)  ,\sigma
^{-1}\left(  k\right)  \right)  $. Therefore, in this case, the number of
inversions of $s_{k}\circ\sigma$ equals the number of inversions of $\sigma$
plus $1$. That is, in this case, $\ell\left(  s_{k}\circ\sigma\right)
=\ell\left(  \sigma\right)  +1$. When $\sigma^{-1}\left(  k\right)
<\sigma^{-1}\left(  k+1\right)  $, a similar argument shows $\ell\left(
s_{k}\circ\sigma\right)  =\ell\left(  \sigma\right)  -1$. This explains why
(\ref{eq.exe.2.5.a.2}) holds (although formalizing this argument will be tedious).

The inversions of $\sigma\circ s_{k}=\left(  4,2,5,1,3\right)  $ are $\left(
1,2\right)  $, $\left(  1,4\right)  $, $\left(  1,5\right)  $, $\left(
2,4\right)  $, $\left(  3,4\right)  $ and $\left(  3,5\right)  $. Unlike the
inversions of $s_{k}\circ\sigma$, these are not directly related to the
inversions of $\sigma$, so the argument in the previous paragraph does not
prove (\ref{eq.exe.2.5.a.1}). However, instead of considering inversions of
$\sigma$, one can consider inversions of $\sigma^{-1}$. These are even more
intuitive: They are the pairs of integers $\left(  i,j\right)  $ with $1\leq
i<j\leq n$ such that $i$ appears further right than $j$ in the one-line
notation of $\sigma$. For instance, the inversions of $\sigma^{-1}$ are
$\left(  1,2\right)  $, $\left(  1,4\right)  $, $\left(  2,4\right)  $,
$\left(  3,4\right)  $ and $\left(  3,5\right)  $, whereas the inversions of
$\left(  \sigma\circ s_{k}\right)  ^{-1}$ are all of these and also $\left(
1,5\right)  $. But there is no need to repeat our proof of
(\ref{eq.exe.2.5.a.2}); it is easier to deduce (\ref{eq.exe.2.5.a.1}) from
(\ref{eq.exe.2.5.a.2}) by applying (\ref{eq.exe.2.5.a.2}) to $\sigma^{-1}$
instead of $\sigma$ and appealing to Exercise \ref{exe.ps2.2.5} \textbf{(f)}.
(Again, see the solution below for the details.)
\end{example}

Notice that Exercise \ref{exe.ps2.2.5} \textbf{(e)} immediately yields
Exercise \ref{exe.ps2.2.4} \textbf{(b)}.

\begin{remark}
When $n=0$ or $n=1$, we have $\left\{  1,2,\ldots,n-1\right\}  =\varnothing$.
Hence, Exercise \ref{exe.ps2.2.4} \textbf{(e)} looks strange in the case when
$n=0$ or $n=1$, because in this case, there are no permutations of the form
$s_{k}$ to begin with. Nevertheless, it is correct. Indeed, when $n=0$ or
$n=1$, there is only one permutation $\sigma\in S_{n}$, namely the identity
permutation $\operatorname*{id}$, and it has length $\ell\left(
\sigma\right)  =\ell\left(  \operatorname*{id}\right)  =0$. Thus, in this
case, Exercise \ref{exe.ps2.2.4} \textbf{(e)} claims that $\operatorname*{id}$
can be written as a composition of $0$ permutations of the form $s_{k}$ (with
$k\in\left\{  1,2,\ldots,n-1\right\}  $). This is true: Even from an empty set
we can always pick $0$ elements; and the composition of $0$ permutations will
be $\operatorname*{id}$.
\end{remark}

\begin{remark}
The word \textquotedblleft length\textquotedblright\ for $\ell\left(
\sigma\right)  $ can be confusing: It does not refer to the length of the
$n$-tuple $\left(  \sigma\left(  1\right)  ,\sigma\left(  2\right)
,\ldots,\sigma\left(  n\right)  \right)  $ (which is $n$). The reason why it
is called \textquotedblleft length\textquotedblright\ is Exercise
\ref{exe.ps2.2.5} \textbf{(g)}: it says that $\ell\left(  \sigma\right)  $ is
the smallest number of permutations of the form $s_{k}$ which can be
multiplied to give $\sigma$; thus, it is the smallest possible length of an
expression of $\sigma$ as a product of $s_{k}$'s.

The use of the word \textquotedblleft length\textquotedblright, unfortunately,
is not standard across literature. Some authors call \textquotedblleft Coxeter
length\textquotedblright\ what we call \textquotedblleft
length\textquotedblright, and use the word \textquotedblleft
length\textquotedblright\ itself for a different notion.
\end{remark}

\begin{exercise}
\label{exe.ps2.2.6}Let $n\in\mathbb{N}$. Let $\sigma\in S_{n}$. In Exercise
\ref{exe.ps2.2.4} \textbf{(b)}, we have seen that $\sigma$ can be written as a
composition of several permutations of the form $s_{k}$ (with $k\in\left\{
1,2,\ldots,n-1\right\}  $). Usually there will be several ways to do so (for
instance, $\operatorname*{id}=s_{1}\circ s_{1}=s_{2}\circ s_{2}=\cdots
=s_{n-1}\circ s_{n-1}$). Show that, whichever of these ways we take, the
number of permutations composed will be congruent to $\ell\left(
\sigma\right)  $ modulo $2$.
\end{exercise}

\subsection{\label{sect.sign}The sign of a permutation}

\begin{definition}
\label{def.perm.sign}Let $n\in\mathbb{N}$.

\textbf{(a)} We define the \textit{sign} of a permutation $\sigma\in S_{n}$ as
the integer $\left(  -1\right)  ^{\ell\left(  \sigma\right)  }$. We denote
this sign by $\left(  -1\right)  ^{\sigma}$ or $\operatorname*{sign}\sigma$ or
$\operatorname*{sgn}\sigma$.

\textbf{(b)} We say that a permutation $\sigma$ is \textit{even} if its sign
is $1$ (that is, if $\ell\left(  \sigma\right)  $ is even), and \textit{odd}
if its sign is $-1$ (that is, if $\ell\left(  \sigma\right)  $ is odd).
\end{definition}

Signs of permutations have the following properties:

\begin{proposition}
\label{prop.perm.signs.basics}Let $n\in\mathbb{N}$.

\textbf{(a)} The sign of the identity permutation $\operatorname*{id}\in
S_{n}$ is $\left(  -1\right)  ^{\operatorname*{id}}=1$. In other words,
$\operatorname*{id}\in S_{n}$ is even.

\textbf{(b)} For every $k\in\left\{  1,2,\ldots,n-1\right\}  $, the sign of
the permutation $s_{k}\in S_{n}$ is $\left(  -1\right)  ^{s_{k}}=-1$.

\textbf{(c)} If $\sigma$ and $\tau$ are two permutations in $S_{n}$, then
$\left(  -1\right)  ^{\sigma\circ\tau}=\left(  -1\right)  ^{\sigma}%
\cdot\left(  -1\right)  ^{\tau}$.

\textbf{(d)} If $\sigma\in S_{n}$, then $\left(  -1\right)  ^{\sigma^{-1}%
}=\left(  -1\right)  ^{\sigma}$. (Here and in the following, the expression
\textquotedblleft$\left(  -1\right)  ^{\sigma^{-1}}$\textquotedblright\ should
be read as \textquotedblleft$\left(  -1\right)  ^{\left(  \sigma^{-1}\right)
}$\textquotedblright, not as \textquotedblleft$\left(  \left(  -1\right)
^{\sigma}\right)  ^{-1}$\textquotedblright; this is similar to Convention
\ref{conv.triple-power} (although $\sigma$ is not a number).)
\end{proposition}

\begin{proof}
[Proof of Proposition \ref{prop.perm.signs.basics}.]\textbf{(a)} The identity
permutation $\operatorname*{id}$ satisfies $\ell\left(  \operatorname*{id}%
\right)  =0$\ \ \ \ \footnote{\textit{Proof.} An inversion of
$\operatorname*{id}$ is the same as a pair $\left(  i,j\right)  $ of integers
satisfying $1\leq i<j\leq n$ and $\operatorname*{id}\left(  i\right)
>\operatorname*{id}\left(  j\right)  $ (by the definition of \textquotedblleft
inversion\textquotedblright). Thus, if $\left(  i,j\right)  $ is an inversion
of $\operatorname*{id}$, then $1\leq i<j\leq n$ and $\operatorname*{id}\left(
i\right)  >\operatorname*{id}\left(  j\right)  $; but this leads to a
contradiction (since $\operatorname*{id}\left(  i\right)  >\operatorname*{id}%
\left(  j\right)  $ contradicts $\operatorname*{id}\left(  i\right)
=i<j=\operatorname*{id}\left(  j\right)  $). Hence, we obtain a contradiction
for each inversion of $\operatorname*{id}$. Thus, there are no inversions of
$\operatorname*{id}$. But $\ell\left(  \operatorname*{id}\right)  $ is defined
as the number of inversions of $\operatorname*{id}$. Hence, $\ell\left(
\operatorname*{id}\right)  =\left(  \text{the number of inversions of
}\operatorname*{id}\right)  =0$ (since there are no inversions of
$\operatorname*{id}$).}. Now, the definition of $\left(  -1\right)
^{\operatorname*{id}}$ yields $\left(  -1\right)  ^{\operatorname*{id}%
}=\left(  -1\right)  ^{\ell\left(  \operatorname*{id}\right)  }=1$ (since
$\ell\left(  \operatorname*{id}\right)  =0$). In other words,
$\operatorname*{id}\in S_{n}$ is even. This proves Proposition
\ref{prop.perm.signs.basics} \textbf{(a)}.

\textbf{(b)} Let $k\in\left\{  1,2,\ldots,n-1\right\}  $. Applying
(\ref{eq.exe.2.5.a.1}) to $\sigma=\operatorname*{id}$, we obtain
\begin{align*}
\ell\left(  \operatorname*{id}\circ s_{k}\right)   &  =%
\begin{cases}
\ell\left(  \operatorname*{id}\right)  +1, & \text{if }\operatorname*{id}%
\left(  k\right)  <\operatorname*{id}\left(  k+1\right)  ;\\
\ell\left(  \operatorname*{id}\right)  -1, & \text{if }\operatorname*{id}%
\left(  k\right)  >\operatorname*{id}\left(  k+1\right)
\end{cases}
\\
&  =\underbrace{\ell\left(  \operatorname*{id}\right)  }_{=0}%
+1\ \ \ \ \ \ \ \ \ \ \left(  \text{since }\operatorname*{id}\left(  k\right)
=k<k+1=\operatorname*{id}\left(  k+1\right)  \right) \\
&  =1.
\end{align*}
This rewrites as $\ell\left(  s_{k}\right)  =1$ (since $\operatorname*{id}%
\circ s_{k}=s_{k}$). Now, the definition of $\left(  -1\right)  ^{s_{k}}$
yields $\left(  -1\right)  ^{s_{k}}=\left(  -1\right)  ^{\ell\left(
s_{k}\right)  }=-1$ (since $\ell\left(  s_{k}\right)  =1$). This proves
Proposition \ref{prop.perm.signs.basics} \textbf{(b)}.

\textbf{(c)} Let $\sigma\in S_{n}$ and $\tau\in S_{n}$. Exercise
\ref{exe.ps2.2.5} \textbf{(b)} yields $\ell\left(  \sigma\circ\tau\right)
\equiv\ell\left(  \sigma\right)  +\ell\left(  \tau\right)  \operatorname{mod}%
2$, so that $\left(  -1\right)  ^{\ell\left(  \sigma\circ\tau\right)
}=\left(  -1\right)  ^{\ell\left(  \sigma\right)  +\ell\left(  \tau\right)
}=\left(  -1\right)  ^{\ell\left(  \sigma\right)  }\cdot\left(  -1\right)
^{\ell\left(  \tau\right)  }$. But the definition of the sign of a permutation
yields $\left(  -1\right)  ^{\sigma\circ\tau}=\left(  -1\right)  ^{\ell\left(
\sigma\circ\tau\right)  }$ and $\left(  -1\right)  ^{\sigma}=\left(
-1\right)  ^{\ell\left(  \sigma\right)  }$ and $\left(  -1\right)  ^{\tau
}=\left(  -1\right)  ^{\ell\left(  \tau\right)  }$. Hence, $\left(  -1\right)
^{\sigma\circ\tau}=\left(  -1\right)  ^{\ell\left(  \sigma\circ\tau\right)
}=\underbrace{\left(  -1\right)  ^{\ell\left(  \sigma\right)  }}_{=\left(
-1\right)  ^{\sigma}}\cdot\underbrace{\left(  -1\right)  ^{\ell\left(
\tau\right)  }}_{=\left(  -1\right)  ^{\tau}}=\left(  -1\right)  ^{\sigma
}\cdot\left(  -1\right)  ^{\tau}$. This proves Proposition
\ref{prop.perm.signs.basics} \textbf{(c)}.

\textbf{(d)} Let $\sigma\in S_{n}$. The definition of $\left(  -1\right)
^{\sigma^{-1}}$ yields $\left(  -1\right)  ^{\sigma^{-1}}=\left(  -1\right)
^{\ell\left(  \sigma^{-1}\right)  }$. But Exercise \ref{exe.ps2.2.5}
\textbf{(f)} says that $\ell\left(  \sigma\right)  =\ell\left(  \sigma
^{-1}\right)  $. The definition of $\left(  -1\right)  ^{\sigma}$ yields
$\left(  -1\right)  ^{\sigma}=\left(  -1\right)  ^{\ell\left(  \sigma\right)
}=\left(  -1\right)  ^{\ell\left(  \sigma^{-1}\right)  }$ (since $\ell\left(
\sigma\right)  =\ell\left(  \sigma^{-1}\right)  $). Compared with $\left(
-1\right)  ^{\sigma^{-1}}=\left(  -1\right)  ^{\ell\left(  \sigma^{-1}\right)
}$, this yields $\left(  -1\right)  ^{\sigma^{-1}}=\left(  -1\right)
^{\sigma}$. This proves Proposition \ref{prop.perm.signs.basics} \textbf{(d)}.
\end{proof}

If you are familiar with some basic concepts of abstract algebra, then you
will immediately notice that parts \textbf{(a)} and \textbf{(c)} of
Proposition \ref{prop.perm.signs.basics} can be summarized as the statement
that \textquotedblleft sign is a group homomorphism from the group $S_{n}$ to
the multiplicative group $\left\{  1,-1\right\}  $\textquotedblright. In this
statement, \textquotedblleft sign\textquotedblright\ means the map from
$S_{n}$ to $\left\{  1,-1\right\}  $ which sends every permutation $\sigma$ to
its sign $\left(  -1\right)  ^{\sigma}=\left(  -1\right)  ^{\ell\left(
\sigma\right)  }$, and the \textquotedblleft multiplicative group $\left\{
1,-1\right\}  $\textquotedblright\ means the group $\left\{  1,-1\right\}  $
whose binary operation is multiplication.

We have defined the sign of a permutation $\sigma\in S_{n}$. More generally,
it is possible to define the sign of a permutation of an arbitrary finite set
$X$, even though the length of such a permutation is not defined!\footnote{How
does it work? If $X$ is a finite set, then we can always find a bijection
$\phi:X\rightarrow\left\{  1,2,\ldots,n\right\}  $ for some $n\in\mathbb{N}$.
(Constructing such a bijection is tantamount to writing down a list of all
elements of $X$, with no duplicates.) Given such a bijection $\phi$, we can
define the sign of any permutation $\sigma$ of $X$ as follows:%
\begin{equation}
\left(  -1\right)  ^{\sigma}=\left(  -1\right)  ^{\phi\circ\sigma\circ
\phi^{-1}}. \label{eq.ps2.S(X).sign.teaser}%
\end{equation}
Here, the right hand side is well-defined because $\phi\circ\sigma\circ
\phi^{-1}$ is a permutation of $\left\{  1,2,\ldots,n\right\}  $. What is not
immediately obvious is that this sign is independent on the choice of $\phi$,
and that it is a group homomorphism to $\left\{  1,-1\right\}  $ (that is, we
have $\left(  -1\right)  ^{\operatorname*{id}}=1$ and $\left(  -1\right)
^{\sigma\circ\tau}=\left(  -1\right)  ^{\sigma}\cdot\left(  -1\right)  ^{\tau
}$). We will prove these facts further below (in Exercise \ref{exe.ps4.2}).}

\begin{exercise}
\label{exe.ps2.2.7}Let $n\geq2$. Show that the number of even permutations in
$S_{n}$ is $n!/2$, and the number of odd permutations in $S_{n}$ is also
$n!/2$.
\end{exercise}

The sign of a permutation is used in the combinatorial definition of the
determinant. Let us briefly show this definition now; we shall return to it
later (in Chapter \ref{chp.det}) to study it in much more detail.

\begin{definition}
\label{def.det.old}Let $n\in\mathbb{N}$. Let $A=\left(  a_{i,j}\right)
_{1\leq i\leq n,\ 1\leq j\leq n}$ be an $n\times n$-matrix (say, with complex
entries, although this does not matter much -- it suffices that the entries
can be added and multiplied and the axioms of associativity, distributivity,
commutativity, unity etc. hold). The \textit{determinant} $\det A$ of $A$ is
defined as%
\begin{equation}
\sum_{\sigma\in S_{n}}\left(  -1\right)  ^{\sigma}a_{1,\sigma\left(  1\right)
}a_{2,\sigma\left(  2\right)  }\cdots a_{n,\sigma\left(  n\right)  }.
\label{eq.det.old}%
\end{equation}

\end{definition}

Let me try to describe the sum (\ref{eq.det.old}) in slightly more visual
terms: The sum (\ref{eq.det.old}) has $n!$ addends, each of which has the form
\textquotedblleft$\left(  -1\right)  ^{\sigma}$ times a
product\textquotedblright. The product has $n$ factors, which are entries of
$A$, and are chosen in such a way that there is exactly one entry taken from
each row and exactly one from each column. Which precise entries are taken
depends on $\sigma$: namely, for each $i$, we take the $\sigma\left(
i\right)  $-th entry from the $i$-th row.

Convince yourself that the classical formulas%
\begin{align*}
\det\left(
\begin{array}
[c]{c}%
a
\end{array}
\right)   &  =a;\\
\det\left(
\begin{array}
[c]{cc}%
a & b\\
c & d
\end{array}
\right)   &  =ad-bc;\\
\det\left(
\begin{array}
[c]{ccc}%
a & b & c\\
d & e & f\\
g & h & i
\end{array}
\right)   &  =aei+bfg+cdh-ahf-bdi-ceg
\end{align*}
are particular cases of (\ref{eq.det.old}). Whenever $n\geq2$, the sum in
(\ref{eq.det.old}) contains precisely $n!/2$ plus signs and $n!/2$ minus signs
(because of Exercise \ref{exe.ps2.2.7}).

Definition \ref{def.det.old} is merely one of several equivalent definitions
of the determinant. You will probably see two of them in an average linear
algebra class. Each of them has its own advantages and drawbacks. Definition
\ref{def.det.old} is the most direct, assuming that one knows about the sign
of a permutation.

\subsection{\label{sect.infperm}Infinite permutations}

(This section is optional; it explores some technical material which is useful
in combinatorics, but is not necessary for what follows. I advise the reader
to skip it at the first read.)

We have introduced the notion of a permutation of an arbitrary set; but so
far, we have only studied permutations of finite sets. In this section (which
is tangential to our project; probably nothing from this section will be used
ever after), let me discuss permutations of the infinite set $\left\{
1,2,3,\ldots\right\}  $. (A lot of what I say below can be easily adapted to
the sets $\mathbb{N}$ and $\mathbb{Z}$ as well.)

We recall that a permutation of a set $X$ means a bijection from $X$ to $X$.

Let $S_{\infty}$ be the symmetric group of the set $\left\{  1,2,3,\ldots
\right\}  $. This is the set of all permutations of $\left\{  1,2,3,\ldots
\right\}  $. It contains the identity permutation $\operatorname*{id}\in
S_{\infty}$ which sends every $i\in\left\{  1,2,3,\ldots\right\}  $ to $i$.
The set $S_{\infty}$ is uncountable\footnote{More generally, while a finite
set of size $n$ has $n!$ permutations, an infinite set $X$ has uncountably
many permutations (even if $X$ is countable).}.

We shall try to study $S_{\infty}$ similarly to how we studied $S_{n}$ for
$n\in\mathbb{N}$. However, we soon will notice that the analogy between
$S_{\infty}$ and $S_{n}$ will break down.\footnote{The uncountability of
$S_{\infty}$ is the first hint that $S_{\infty}$ is \textquotedblleft too
large\textquotedblright\ a set to be a good analogue of the finite set $S_{n}%
$.} To amend this, we shall define a subset $S_{\left(  \infty\right)  }$ of
$S_{\infty}$ (mind the parentheses around the \textquotedblleft$\infty
$\textquotedblright) which is smaller and more wieldy, and indeed shares many
of the properties of the finite symmetric group $S_{n}$.

We define $S_{\left(  \infty\right)  }$ as follows:%
\begin{equation}
S_{\left(  \infty\right)  }=\left\{  \sigma\in S_{\infty}\ \mid\ \sigma\left(
i\right)  =i\text{ for all but finitely many }i\in\left\{  1,2,3,\ldots
\right\}  \right\}  . \label{eq.S(infty).def}%
\end{equation}
Let us first explain what \textquotedblleft all but finitely many
$i\in\left\{  1,2,3,\ldots\right\}  $\textquotedblright\ means:

\begin{definition}
\label{def.allbutfin}Let $I$ be a set. Let $\mathcal{A}\left(  i\right)  $ be
a statement for every $i\in I$. Then, we say that \textquotedblleft%
$\mathcal{A}\left(  i\right)  $ for all but finitely many $i\in I$%
\textquotedblright\ if and only if there exists some finite subset $J$ of $I$
such that every $i\in I\setminus J$ satisfies $\mathcal{A}\left(  i\right)  $.\ \ \ \ \footnotemark
\end{definition}

\footnotetext{Thus, the statement \textquotedblleft$\mathcal{A}\left(
i\right)  $ for all but finitely many $i\in I$\textquotedblright\ can be
restated as \textquotedblleft$\mathcal{A}\left(  i\right)  $ holds for all
$i\in I$, apart from finitely many exceptions\textquotedblright\ or as
\textquotedblleft there are only finitely many $i\in I$ which do not satisfy
$\mathcal{A}\left(  i\right)  $\textquotedblright. I prefer the first wording,
because it makes the most sense in constructive logic.
\par
\textbf{Caution:} Do not confuse the words \textquotedblleft all but finitely
many $i\in I$\textquotedblright\ in this definition with the words
\textquotedblleft infinitely many $i\in I$\textquotedblright. For instance, it
is true that $n$ is even for infinitely many $n\in\mathbb{Z}$, but it is not
true that $n$ is even for all but finitely many $n\in\mathbb{Z}$. Conversely,
it is true that $n>1$ for all but finitely many $n\in\left\{  1,2\right\}  $
(because the only $n\in\left\{  1,2\right\}  $ which does not satisfy $n>1$ is
$1$), but it is not true that $n>1$ for infinitely many $n\in\left\{
1,2\right\}  $ (because there are no infinitely many $n\in\left\{
1,2\right\}  $ to begin with).
\par
You will encounter the \textquotedblleft all but finitely
many\textquotedblright\ formulation often in abstract algebra. (Some people
abbreviate it as \textquotedblleft almost all\textquotedblright, but this
abbreviation means other things as well.)} Thus, for a permutation $\sigma\in
S_{\infty}$, we have the following equivalence of statements:%
\begin{align*}
&  \ \left(  \sigma\left(  i\right)  =i\text{ for all but finitely many }%
i\in\left\{  1,2,3,\ldots\right\}  \right) \\
&  \Longleftrightarrow\ \left(  \text{there exists some finite subset }J\text{
of }\left\{  1,2,3,\ldots\right\}  \text{ such that}\right. \\
&  \ \ \ \ \ \ \ \ \ \ \left.  \text{every }i\in\left\{  1,2,3,\ldots\right\}
\setminus J\text{ satisfies }\sigma\left(  i\right)  =i\right) \\
&  \Longleftrightarrow\ \left(  \text{there exists some finite subset }J\text{
of }\left\{  1,2,3,\ldots\right\}  \text{ such that}\right. \\
&  \ \ \ \ \ \ \ \ \ \ \left.  \text{the only }i\in\left\{  1,2,3,\ldots
\right\}  \text{ that satisfy }\sigma\left(  i\right)  \neq i\text{ are
elements of }J\right) \\
&  \Longleftrightarrow\ \left(  \text{the set of all }i\in\left\{
1,2,3,\ldots\right\}  \text{ that satisfy }\sigma\left(  i\right)  \neq
i\text{ is}\right. \\
&  \ \ \ \ \ \ \ \ \ \ \left.  \text{contained in some finite subset }J\text{
of }\left\{  1,2,3,\ldots\right\}  \right) \\
&  \Longleftrightarrow\ \left(  \text{there are only finitely many }%
i\in\left\{  1,2,3,\ldots\right\}  \text{ that satisfy }\sigma\left(
i\right)  \neq i\right)  .
\end{align*}
Hence, (\ref{eq.S(infty).def}) rewrites as follows:%
\[
S_{\left(  \infty\right)  }=\left\{  \sigma\in S_{\infty}\ \mid\ \text{there
are only finitely many }i\in\left\{  1,2,3,\ldots\right\}  \text{ that satisfy
}\sigma\left(  i\right)  \neq i\right\}  .
\]

\begin{example}
Here is an example of a permutation which is in $S_{\infty}$ but not in
$S_{\left(  \infty\right)  }$: Let $\tau$ be the permutation of $\left\{
1,2,3,\ldots\right\}  $ given by%
\[
\left(  \tau\left(  1\right)  ,\tau\left(  2\right)  ,\tau\left(  3\right)
,\tau\left(  4\right)  ,\tau\left(  5\right)  ,\tau\left(  6\right)
,\ldots\right)  =\left(  2,1,4,3,6,5,\ldots\right)  .
\]
(It adds $1$ to every odd positive integer, and subtracts $1$ from every even
positive integer.) Then, $\tau\in S_{\infty}$ but $\tau\notin S_{\left(
\infty\right)  }$.
\end{example}

On the other hand, let us show some examples of permutations in $S_{\left(
\infty\right)  }$. For each $i\in\left\{  1,2,3,\ldots\right\}  $, let $s_{i}$
be the permutation in $S_{\infty}$ that swaps $i$ with $i+1$ but leaves all
other numbers unchanged. (This is similar to the permutation $s_{i}$ in
$S_{n}$ that was defined earlier. We have taken the liberty to re-use the name
$s_{i}$, hoping that no confusion will arise.)

Again, we have $s_{i}^{2}=\operatorname*{id}$ for every $i\in\left\{
1,2,3,\ldots\right\}  $ (where $\alpha^{2}$ means $\alpha\circ\alpha$ for any
$\alpha\in S_{\infty}$).

\begin{proposition}
\label{prop.S(infty).si}We have $s_{k}\in S_{\left(  \infty\right)  }$ for
every $k\in\left\{  1,2,3,\ldots\right\}  $.
\end{proposition}

\begin{proof}
[Proof of Proposition \ref{prop.S(infty).si}.]Let $k\in\left\{  1,2,3,\ldots
\right\}  $. The permutation $s_{k}$ has been defined as the permutation in
$S_{\infty}$ that swaps $k$ with $k+1$ but leaves all other numbers unchanged.
In other words, it satisfies $s_{k}\left(  k\right)  =k+1$, $s_{k}\left(
k+1\right)  =k$ and%
\begin{equation}
s_{k}\left(  i\right)  =i\ \ \ \ \ \ \ \ \ \ \text{for every }i\in\left\{
1,2,3,\ldots\right\}  \text{ such that }i\notin\left\{  k,k+1\right\}  .
\label{pf.prop.S(infty).si.1}%
\end{equation}

Now, every $i\in\left\{  1,2,3,\ldots\right\}  \setminus\left\{
k,k+1\right\}  $ satisfies $s_{k}\left(  i\right)  =i$%
\ \ \ \ \footnote{\textit{Proof.} Let $i\in\left\{  1,2,3,\ldots\right\}
\setminus\left\{  k,k+1\right\}  $. Thus, $i\in\left\{  1,2,3,\ldots\right\}
$ and $i\notin\left\{  k,k+1\right\}  $. Hence, (\ref{pf.prop.S(infty).si.1})
shows that $s_{k}\left(  i\right)  =i$, qed.}. Hence, there exists some finite
subset $J$ of $\left\{  1,2,3,\ldots\right\}  $ such that every $i\in\left\{
1,2,3,\ldots\right\}  \setminus J$ satisfies $s_{k}\left(  i\right)  =i$
(namely, $J=\left\{  k,k+1\right\}  $). In other words, $s_{k}\left(
i\right)  =i$ for all but finitely many $i\in\left\{  1,2,3,\ldots\right\}  $.

Thus, $s_{k}$ is an element of $S_{\infty}$ satisfying $s_{k}\left(  i\right)
=i$ for all but finitely many $i\in\left\{  1,2,3,\ldots\right\}  $. Hence,%
\[
s_{k}\in\left\{  \sigma\in S_{\infty}\ \mid\ \sigma\left(  i\right)  =i\text{
for all but finitely many }i\in\left\{  1,2,3,\ldots\right\}  \right\}
=S_{\left(  \infty\right)  }.
\]
This proves Proposition \ref{prop.S(infty).si}.
\end{proof}

Permutations can be composed and inverted, leading to new permutations. Let us
first see that the same is true for elements of $S_{\left(  \infty\right)  }$:

\begin{proposition}
\label{prop.S(infty).group}\textbf{(a)} The identity permutation
$\operatorname*{id}\in S_{\infty}$ of $\left\{  1,2,3,\ldots\right\}  $
satisfies $\operatorname*{id}\in S_{\left(  \infty\right)  }$.

\textbf{(b)} For every $\sigma\in S_{\left(  \infty\right)  }$ and $\tau\in
S_{\left(  \infty\right)  }$, we have $\sigma\circ\tau\in S_{\left(
\infty\right)  }$.

\textbf{(c)} For every $\sigma\in S_{\left(  \infty\right)  }$, we have
$\sigma^{-1}\in S_{\left(  \infty\right)  }$.
\end{proposition}

\begin{vershort}
\begin{proof}
[Proof of Proposition \ref{prop.S(infty).group}.]We have defined $S_{\left(
\infty\right)  }$ as the set of all $\sigma\in S_{\infty}$ such that
$\sigma\left(  i\right)  =i$ for all but finitely many $i\in\left\{
1,2,3,\ldots\right\}  $. In other words, $S_{\left(  \infty\right)  }$ is the
set of all $\sigma\in S_{\infty}$ such that there exists a finite subset $K$
of $\left\{  1,2,3,\ldots\right\}  $ such that $\left(  \text{every }%
i\in\left\{  1,2,3,\ldots\right\}  \setminus K\text{ satisfies }\sigma\left(
i\right)  =i\right)  $. As a consequence, we have the following two facts:

\begin{itemize}
\item If $K$ is a finite subset of $\left\{  1,2,3,\ldots\right\}  $, and if
$\gamma\in S_{\infty}$ is a permutation such that%
\begin{equation}
\left(  \text{every }i\in\left\{  1,2,3,\ldots\right\}  \setminus K\text{
satisfies }\gamma\left(  i\right)  =i\right)  ,
\label{pf.prop.S(infty).group.short.lem.hyp}%
\end{equation}
then%
\begin{equation}
\gamma\in S_{\left(  \infty\right)  }.
\label{pf.prop.S(infty).group.short.lem}%
\end{equation}

\item If $\gamma\in S_{\left(  \infty\right)  }$, then
\begin{equation}
\left(
\begin{array}
[c]{c}%
\text{there exists some finite subset }K\text{ of }\left\{  1,2,3,\ldots
\right\} \\
\text{such that every }i\in\left\{  1,2,3,\ldots\right\}  \setminus K\text{
satisfies }\gamma\left(  i\right)  =i
\end{array}
\right)  . \label{pf.prop.S(infty).group.short.lem2}%
\end{equation}

\end{itemize}

We can now step to the actual proof of Proposition \ref{prop.S(infty).group}.

\textbf{(a)} Every $i\in\left\{  1,2,3,\ldots\right\}  \setminus\varnothing$
satisfies $\operatorname*{id}\left(  i\right)  =i$. Thus,
(\ref{pf.prop.S(infty).group.short.lem}) (applied to $K=\varnothing$ and
$\gamma=\operatorname*{id}$) yields $\operatorname*{id}\in S_{\left(
\infty\right)  }$. This proves Proposition \ref{prop.S(infty).group}
\textbf{(a)}.

\textbf{(b)} Let $\sigma\in S_{\left(  \infty\right)  }$ and $\tau\in
S_{\left(  \infty\right)  }$.

From (\ref{pf.prop.S(infty).group.short.lem2}) (applied to $\gamma=\sigma$),
we conclude that there exists some finite subset $K$ of $\left\{
1,2,3,\ldots\right\}  $ such that every $i\in\left\{  1,2,3,\ldots\right\}
\setminus K$ satisfies $\sigma\left(  i\right)  =i$. Let us denote this $K$ by
$J_{1}$. Thus, $J_{1}$ is a finite subset of $\left\{  1,2,3,\ldots\right\}
$, and%
\begin{equation}
\text{every }i\in\left\{  1,2,3,\ldots\right\}  \setminus J_{1}\text{
satisfies }\sigma\left(  i\right)  =i.
\label{pf.prop.S(infty).group.short.b.1}%
\end{equation}

From (\ref{pf.prop.S(infty).group.short.lem2}) (applied to $\gamma=\tau$), we
conclude that there exists some finite subset $K$ of $\left\{  1,2,3,\ldots
\right\}  $ such that every $i\in\left\{  1,2,3,\ldots\right\}  \setminus K$
satisfies $\tau\left(  i\right)  =i$. Let us denote this $K$ by $J_{2}$. Thus,
$J_{2}$ is a finite subset of $\left\{  1,2,3,\ldots\right\}  $, and%
\begin{equation}
\text{every }i\in\left\{  1,2,3,\ldots\right\}  \setminus J_{2}\text{
satisfies }\tau\left(  i\right)  =i. \label{pf.prop.S(infty).group.short.b.2}%
\end{equation}

The sets $J_{1}$ and $J_{2}$ are finite. Hence, their union $J_{1}\cup J_{2}$
is finite. Moreover,
\[
\text{every }i\in\left\{  1,2,3,\ldots\right\}  \setminus\left(  J_{1}\cup
J_{2}\right)  \text{ satisfies }\left(  \sigma\circ\tau\right)  \left(
i\right)  =i
\]
\footnote{\textit{Proof.} Let $i\in\left\{  1,2,3,\ldots\right\}
\setminus\left(  J_{1}\cup J_{2}\right)  $. Thus, $i\in\left\{  1,2,3,\ldots
\right\}  $ and $i\notin J_{1}\cup J_{2}$.
\par
We have $i\notin J_{1}\cup J_{2}$ and thus $i\notin J_{1}$ (since
$J_{1}\subseteq J_{1}\cup J_{2}$). Hence, $i\in\left\{  1,2,3,\ldots\right\}
\setminus J_{1}$. Similarly, $i\in\left\{  1,2,3,\ldots\right\}  \setminus
J_{2}$. Thus, (\ref{pf.prop.S(infty).group.short.b.2}) yields $\tau\left(
i\right)  =i$. Hence, $\left(  \sigma\circ\tau\right)  \left(  i\right)
=\sigma\left(  \underbrace{\tau\left(  i\right)  }_{=i}\right)  =\sigma\left(
i\right)  =i$ (by (\ref{pf.prop.S(infty).group.short.b.1})), qed.}. Therefore,
(\ref{pf.prop.S(infty).group.short.lem}) (applied to $K=J_{1}\cup J_{2}$ and
$\gamma=\sigma\circ\tau$) yields $\sigma\circ\tau\in S_{\left(  \infty\right)
}$. This proves Proposition \ref{prop.S(infty).group} \textbf{(b)}.

\textbf{(c)} Let $\sigma\in S_{\left(  \infty\right)  }$.

From (\ref{pf.prop.S(infty).group.short.lem2}) (applied to $\gamma=\sigma$),
we conclude that there exists some finite subset $K$ of $\left\{
1,2,3,\ldots\right\}  $ such that every $i\in\left\{  1,2,3,\ldots\right\}
\setminus K$ satisfies $\sigma\left(  i\right)  =i$. Consider this $K$. Thus,
$K$ is a finite subset of $\left\{  1,2,3,\ldots\right\}  $, and%
\begin{equation}
\text{every }i\in\left\{  1,2,3,\ldots\right\}  \setminus K\text{ satisfies
}\sigma\left(  i\right)  =i. \label{pf.prop.S(infty).group.short.c.3}%
\end{equation}

Now,
\[
\text{every }i\in\left\{  1,2,3,\ldots\right\}  \setminus K\text{ satisfies
}\sigma^{-1}\left(  i\right)  =i
\]
\footnote{\textit{Proof.} Let $i\in\left\{  1,2,3,\ldots\right\}  \setminus
K$. Thus, $\sigma\left(  i\right)  =i$ (according to
(\ref{pf.prop.S(infty).group.short.c.3})), so that $\sigma^{-1}\left(
i\right)  =i$, qed.}. Therefore, (\ref{pf.prop.S(infty).group.short.lem})
(applied to $\gamma=\sigma^{-1}$) yields $\sigma^{-1}\in S_{\left(
\infty\right)  }$. This proves Proposition \ref{prop.S(infty).group}
\textbf{(c)}.
\end{proof}
\end{vershort}

\begin{verlong}
\begin{proof}
[Proof of Proposition \ref{prop.S(infty).group}.]Let us record the following facts:

\begin{itemize}
\item If $K$ is a finite subset of $\left\{  1,2,3,\ldots\right\}  $, and if
$\gamma\in S_{\infty}$ is a permutation such that%
\begin{equation}
\left(  \text{every }i\in\left\{  1,2,3,\ldots\right\}  \setminus K\text{
satisfies }\gamma\left(  i\right)  =i\right)  ,
\label{pf.prop.S(infty).group.lem.hyp}%
\end{equation}
then%
\begin{equation}
\gamma\in S_{\left(  \infty\right)  }. \label{pf.prop.S(infty).group.lem}%
\end{equation}
\footnote{\textit{Proof of (\ref{pf.prop.S(infty).group.lem}):} Let $K$ be a
finite subset of $\left\{  1,2,3,\ldots\right\}  $, and let $\gamma\in
S_{\infty}$ be a permutation such that (\ref{pf.prop.S(infty).group.lem.hyp})
holds. We need to prove that $\gamma\in S_{\left(  \infty\right)  }$.
\par
Every $i\in\left\{  1,2,3,\ldots\right\}  \setminus K$ satisfies
$\gamma\left(  i\right)  =i$ (since (\ref{pf.prop.S(infty).group.lem.hyp})
holds). Thus, there exists some finite subset $J$ of $\left\{  1,2,3,\ldots
\right\}  $ such that every $i\in\left\{  1,2,3,\ldots\right\}  \setminus J$
satisfies $\gamma\left(  i\right)  =i$ (namely, $J=K$). In other words,
$\gamma\left(  i\right)  =i$ for all but finitely many $i\in\left\{
1,2,3,\ldots\right\}  $.
\par
Thus, $\gamma$ is an element of $S_{\infty}$ satisfying $\gamma\left(
i\right)  =i$ for all but finitely many $i\in\left\{  1,2,3,\ldots\right\}  $.
Hence,%
\[
\gamma\in\left\{  \sigma\in S_{\infty}\ \mid\ \sigma\left(  i\right)  =i\text{
for all but finitely many }i\in\left\{  1,2,3,\ldots\right\}  \right\}
=S_{\left(  \infty\right)  }.
\]
This proves (\ref{pf.prop.S(infty).group.lem}).}

\item If $\gamma\in S_{\left(  \infty\right)  }$, then
\begin{equation}
\left(
\begin{array}
[c]{c}%
\text{there exists some finite subset }J\text{ of }\left\{  1,2,3,\ldots
\right\} \\
\text{such that every }i\in\left\{  1,2,3,\ldots\right\}  \setminus J\text{
satisfies }\gamma\left(  i\right)  =i
\end{array}
\right)  \label{pf.prop.S(infty).group.lem2}%
\end{equation}
\footnote{\textit{Proof of (\ref{pf.prop.S(infty).group.lem2}):} Let
$\gamma\in S_{\left(  \infty\right)  }$. Then,
\[
\gamma\in S_{\left(  \infty\right)  }=\left\{  \sigma\in S_{\infty}%
\ \mid\ \sigma\left(  i\right)  =i\text{ for all but finitely many }%
i\in\left\{  1,2,3,\ldots\right\}  \right\}  .
\]
In other words, $\gamma$ is an element of $S_{\infty}$ such that
$\gamma\left(  i\right)  =i$ for all but finitely many $i\in\left\{
1,2,3,\ldots\right\}  $.
\par
Now, we know that $\gamma\left(  i\right)  =i$ for all but finitely many
$i\in\left\{  1,2,3,\ldots\right\}  $. In other words, there exists some
finite subset $J$ of $\left\{  1,2,3,\ldots\right\}  $ such that every
$i\in\left\{  1,2,3,\ldots\right\}  \setminus J$ satisfies $\gamma\left(
i\right)  =i$. This proves (\ref{pf.prop.S(infty).group.lem2}).}.
\end{itemize}

We can now step to the actual proof of Proposition \ref{prop.S(infty).group}.

\textbf{(a)} Every $i\in\left\{  1,2,3,\ldots\right\}  \setminus\varnothing$
satisfies $\operatorname*{id}\left(  i\right)  =i$. Thus,
(\ref{pf.prop.S(infty).group.lem}) (applied to $K=\varnothing$ and
$\gamma=\operatorname*{id}$) yields $\operatorname*{id}\in S_{\left(
\infty\right)  }$. This proves Proposition \ref{prop.S(infty).group}
\textbf{(a)}.

\textbf{(b)} Let $\sigma\in S_{\left(  \infty\right)  }$ and $\tau\in
S_{\left(  \infty\right)  }$.

From (\ref{pf.prop.S(infty).group.lem2}) (applied to $\gamma=\sigma$), we
conclude that there exists some finite subset $J$ of $\left\{  1,2,3,\ldots
\right\}  $ such that every $i\in\left\{  1,2,3,\ldots\right\}  \setminus J$
satisfies $\sigma\left(  i\right)  =i$. Let us denote this $J$ by $J_{1}$.
Thus, $J_{1}$ is a finite subset of $\left\{  1,2,3,\ldots\right\}  $, and%
\begin{equation}
\text{every }i\in\left\{  1,2,3,\ldots\right\}  \setminus J_{1}\text{
satisfies }\sigma\left(  i\right)  =i. \label{pf.prop.S(infty).group.b.1}%
\end{equation}

From (\ref{pf.prop.S(infty).group.lem2}) (applied to $\gamma=\tau$), we
conclude that there exists some finite subset $J$ of $\left\{  1,2,3,\ldots
\right\}  $ such that every $i\in\left\{  1,2,3,\ldots\right\}  \setminus J$
satisfies $\tau\left(  i\right)  =i$. Let us denote this $J$ by $J_{2}$. Thus,
$J_{2}$ is a finite subset of $\left\{  1,2,3,\ldots\right\}  $, and%
\begin{equation}
\text{every }i\in\left\{  1,2,3,\ldots\right\}  \setminus J_{2}\text{
satisfies }\tau\left(  i\right)  =i. \label{pf.prop.S(infty).group.b.2}%
\end{equation}

The sets $J_{1}$ and $J_{2}$ are finite. Hence, their union $J_{1}\cup J_{2}$
is finite. Moreover,
\[
\text{every }i\in\left\{  1,2,3,\ldots\right\}  \setminus\left(  J_{1}\cup
J_{2}\right)  \text{ satisfies }\left(  \sigma\circ\tau\right)  \left(
i\right)  =i
\]
\footnote{\textit{Proof.} Let $i\in\left\{  1,2,3,\ldots\right\}
\setminus\left(  J_{1}\cup J_{2}\right)  $. Thus, $i\in\left\{  1,2,3,\ldots
\right\}  $ and $i\notin J_{1}\cup J_{2}$.
\par
If we had $i\in J_{1}$, then we would have $i\in J_{1}\subseteq J_{1}\cup
J_{2}$, which would contradict $i\notin J_{1}\cup J_{2}$. Thus, we cannot have
$i\in J_{1}$. In other words, we have $i\notin J_{1}$. Hence, $i\in\left\{
1,2,3,\ldots\right\}  \setminus J_{1}$. Similarly, $i\in\left\{
1,2,3,\ldots\right\}  \setminus J_{2}$. Thus,
(\ref{pf.prop.S(infty).group.b.2}) yields $\tau\left(  i\right)  =i$. Hence,
$\left(  \sigma\circ\tau\right)  \left(  i\right)  =\sigma\left(
\underbrace{\tau\left(  i\right)  }_{=i}\right)  =\sigma\left(  i\right)  =i$
(by (\ref{pf.prop.S(infty).group.b.1})), qed.}. Therefore,
(\ref{pf.prop.S(infty).group.lem}) (applied to $K=J_{1}\cup J_{2}$ and
$\gamma=\sigma\circ\tau$) yields $\sigma\circ\tau\in S_{\left(  \infty\right)
}$. This proves Proposition \ref{prop.S(infty).group} \textbf{(b)}.

\textbf{(c)} Let $\sigma\in S_{\left(  \infty\right)  }$.

From (\ref{pf.prop.S(infty).group.lem2}) (applied to $\gamma=\sigma$), we
conclude that there exists some finite subset $J$ of $\left\{  1,2,3,\ldots
\right\}  $ such that every $i\in\left\{  1,2,3,\ldots\right\}  \setminus J$
satisfies $\sigma\left(  i\right)  =i$. Consider this $J$. Thus, $J$ is a
finite subset of $\left\{  1,2,3,\ldots\right\}  $, and%
\begin{equation}
\text{every }i\in\left\{  1,2,3,\ldots\right\}  \setminus J\text{ satisfies
}\sigma\left(  i\right)  =i. \label{pf.prop.S(infty).group.c.3}%
\end{equation}

Now,
\[
\text{every }i\in\left\{  1,2,3,\ldots\right\}  \setminus J\text{ satisfies
}\sigma^{-1}\left(  i\right)  =i
\]
\footnote{\textit{Proof.} Let $i\in\left\{  1,2,3,\ldots\right\}  \setminus
J$. Thus, $\sigma\left(  i\right)  =i$ (according to
(\ref{pf.prop.S(infty).group.c.3})), so that $\sigma^{-1}\left(  i\right)
=i$, qed.}. Therefore, (\ref{pf.prop.S(infty).group.lem}) (applied to $K=J$
and $\gamma=\sigma^{-1}$) yields $\sigma^{-1}\in S_{\left(  \infty\right)  }$.
This proves Proposition \ref{prop.S(infty).group} \textbf{(c)}.
\end{proof}
\end{verlong}

In the language of group theorists, Proposition \ref{prop.S(infty).group} show
that $S_{\left(  \infty\right)  }$ is a subgroup of the group $S_{\infty}$.
The elements of $S_{\left(  \infty\right)  }$ are called the \textit{finitary
permutations of }$\left\{  1,2,3,\ldots\right\}  $, and $S_{\left(
\infty\right)  }$ is called the \textit{finitary symmetric group of }$\left\{
1,2,3,\ldots\right\}  $.

We now have the following analogue of Exercise \ref{exe.ps2.2.4} (without its
part \textbf{(c)}):

\begin{exercise}
\label{exe.ps2.2.4'}\textbf{(a)} Show that $s_{i}\circ s_{i+1}\circ
s_{i}=s_{i+1}\circ s_{i}\circ s_{i+1}$ for all $i\in\left\{  1,2,3,\ldots
\right\}  $.

\textbf{(b)} Show that every permutation $\sigma\in S_{\left(  \infty\right)
}$ can be written as a composition of several permutations of the form $s_{k}$
(with $k\in\left\{  1,2,3,\ldots\right\}  $).
\end{exercise}

\begin{remark}
In the language of group theory, the statement of Exercise \ref{exe.ps2.2.4'}
\textbf{(b)} says (or, more precisely, yields) that the permutations
$s_{1},s_{2},s_{3},\ldots$ generate the group $S_{\left(  \infty\right)  }$.
\end{remark}

If $\sigma\in S_{\infty}$ is a permutation, then an \textit{inversion} of
$\sigma$ means a pair $\left(  i,j\right)  $ of integers satisfying $1\leq
i<j$ and $\sigma\left(  i\right)  >\sigma\left(  j\right)  $. This definition
of an inversion is, of course, analogous to the definition of an inversion of
a $\sigma\in S_{n}$; thus we could try to define the length of a $\sigma\in
S_{\infty}$. However, here we run into troubles: A permutation $\sigma\in
S_{\infty}$ might have infinitely many inversions!

It is here that we really need to restrict ourselves to $S_{\left(
\infty\right)  }$. This indeed helps:

\begin{proposition}
\label{prop.S(infty).inversions}Let $\sigma\in S_{\left(  \infty\right)  }$. Then:

\textbf{(a)} There exists some $N\in\left\{  1,2,3,\ldots\right\}  $ such that
every integer $i>N$ satisfies $\sigma\left(  i\right)  =i$.

\textbf{(b)} There are only finitely many inversions of $\sigma$.
\end{proposition}

\begin{vershort}
\begin{proof}
[Proof of Proposition \ref{prop.S(infty).inversions}.]\textbf{(a)} We can
apply (\ref{pf.prop.S(infty).group.short.lem2}) to $\gamma=\sigma$. As a
consequence, we obtain that there exists some finite subset $K$ of $\left\{
1,2,3,\ldots\right\}  $ such that%
\begin{equation}
\text{every }i\in\left\{  1,2,3,\ldots\right\}  \setminus K\text{ satisfies
}\sigma\left(  i\right)  =i. \label{pf.prop.S(infty).inversions.short.a.J}%
\end{equation}
Consider this $K$.

The set $K$ is finite. Hence, the set $K\cup\left\{  1\right\}  $ is finite;
this set is also nonempty (since it contains $1$) and a subset of $\left\{
1,2,3,\ldots\right\}  $. Therefore, this set $K\cup\left\{  1\right\}  $ has a
greatest element (since every finite nonempty subset of $\left\{
1,2,3,\ldots\right\}  $ has a greatest element). Let $n$ be this greatest
element. Clearly, $n\in K\cup\left\{  1\right\}  \subseteq\left\{
1,2,3,\ldots\right\}  $, so that $n>0$.

Every $j\in K\cup\left\{  1\right\}  $ satisfies
\begin{equation}
j\leq n \label{pf.prop.S(infty).inversions.short.a.1}%
\end{equation}
(since $n$ is the greatest element of $K\cup\left\{  1\right\}  $). Now, let
$i$ be an integer such that $i>n$. Then, $i>n>0$, so that $i$ is a positive
integer. If we had $i\in K$, then we would have $i\in K\subseteq K\cup\left\{
1\right\}  $ and thus $i\leq n$ (by
(\ref{pf.prop.S(infty).inversions.short.a.1}), applied to $j=i$), which would
contradict $i>n$. Hence, we cannot have $i\in K$. We thus have $i\notin K$.
Since $i\in\left\{  1,2,3,\ldots\right\}  $, this shows that $i\in\left\{
1,2,3,\ldots\right\}  \setminus K$. Thus, $\sigma\left(  i\right)  =i$ (by
(\ref{pf.prop.S(infty).inversions.short.a.J})).

Let us now forget that we fixed $i$. We thus have shown that every integer
$i>n$ satisfies $\sigma\left(  i\right)  =i$. Hence, Proposition
\ref{prop.S(infty).inversions} \textbf{(a)} holds (we can take $N=n$).

\textbf{(b)} Proposition \ref{prop.S(infty).inversions} \textbf{(a)} shows
that there exists some $N\in\left\{  1,2,3,\ldots\right\}  $ such that%
\begin{equation}
\text{every integer }i>N\text{ satisfies }\sigma\left(  i\right)  =i.
\label{pf.prop.S(infty).inversions.short.b.N}%
\end{equation}
Consider such an $N$. We shall now show that
\[
\text{every inversion of }\sigma\text{ is an element of }\left\{
1,2,\ldots,N\right\}  ^{2}.
\]

In fact, let $c$ be an inversion of $\sigma$. We will show that $c$ is an
element of $\left\{  1,2,\ldots,N\right\}  ^{2}$.

We know that $c$ is an inversion of $\sigma$. In other words, $c$ is a pair
$\left(  i,j\right)  $ of integers satisfying $1\leq i<j$ and $\sigma\left(
i\right)  >\sigma\left(  j\right)  $ (by the definition of an
\textquotedblleft inversion of $\sigma$\textquotedblright). Consider this
$\left(  i,j\right)  $. We then have $i\leq N$\ \ \ \footnote{\textit{Proof.}
Assume the contrary. Thus, $i>N$. Hence,
(\ref{pf.prop.S(infty).inversions.short.b.N}) shows that $\sigma\left(
i\right)  =i$. Also, $i<j$, so that $j>i>N$. Hence,
(\ref{pf.prop.S(infty).inversions.short.b.N}) (applied to $j$ instead of $i$)
shows that $\sigma\left(  j\right)  =j$. Thus, $\sigma\left(  i\right)
=i<j=\sigma\left(  j\right)  $. This contradicts $\sigma\left(  i\right)
>\sigma\left(  j\right)  $. This contradiction shows that our assumption was
wrong, qed.} and $j\leq N$\ \ \ \ \footnote{\textit{Proof.} Assume the
contrary. Thus, $j>N$. Hence, (\ref{pf.prop.S(infty).inversions.short.b.N})
(applied to $j$ instead of $i$) shows that $\sigma\left(  j\right)  =j$. Now,
$\sigma\left(  i\right)  >\sigma\left(  j\right)  =j>N$. Therefore,
(\ref{pf.prop.S(infty).inversions.short.b.N}) (applied to $\sigma\left(
i\right)  $ instead of $i$) yields $\sigma\left(  \sigma\left(  i\right)
\right)  =\sigma\left(  i\right)  $. But $\sigma$ is a permutation, and thus
an injective map. Hence, from $\sigma\left(  \sigma\left(  i\right)  \right)
=\sigma\left(  i\right)  $, we obtain $\sigma\left(  i\right)  =i$. Thus,
$\sigma\left(  i\right)  =i<j=\sigma\left(  j\right)  $. This contradicts
$\sigma\left(  i\right)  >\sigma\left(  j\right)  $. This contradiction shows
that our assumption was wrong, qed.}. Consequently, $\left(  i,j\right)
\in\left\{  1,2,\ldots,N\right\}  ^{2}$. Hence, $c=\left(  i,j\right)
\in\left\{  1,2,\ldots,N\right\}  ^{2}$.

Now, let us forget that we fixed $c$. We thus have shown that if $c$ is an
inversion of $\sigma$, then $c$ is an element of $\left\{  1,2,\ldots
,N\right\}  ^{2}$. In other words, every inversion of $\sigma$ is an element
of $\left\{  1,2,\ldots,N\right\}  ^{2}$. Thus, there are only finitely many
inversions of $\sigma$ (since there are only finitely many elements of
$\left\{  1,2,\ldots,N\right\}  ^{2}$). Proposition
\ref{prop.S(infty).inversions} \textbf{(b)} is thus proven.
\end{proof}
\end{vershort}

\begin{verlong}
\begin{proof}
[Proof of Proposition \ref{prop.S(infty).inversions}.]\textbf{(a)} We can
apply (\ref{pf.prop.S(infty).group.lem2}) to $\gamma=\sigma$. As a
consequence, we obtain that there exists some finite subset $J$ of $\left\{
1,2,3,\ldots\right\}  $ such that%
\begin{equation}
\text{every }i\in\left\{  1,2,3,\ldots\right\}  \setminus J\text{ satisfies
}\sigma\left(  i\right)  =i. \label{pf.prop.S(infty).inversions.a.J}%
\end{equation}
Consider this $J$.

The set $J$ is finite. Hence, the set $J\cup\left\{  1\right\}  $ is finite;
this set is also nonempty (since it contains $1$) and a subset of $\left\{
1,2,3,\ldots\right\}  $. Therefore, this set $J\cup\left\{  1\right\}  $ has a
greatest element (since every finite nonempty subset of $\left\{
1,2,3,\ldots\right\}  $ has a greatest element). Let $n$ be this greatest
element. Clearly, $n\in J\cup\left\{  1\right\}  \subseteq\left\{
1,2,3,\ldots\right\}  $, so that $n>0$.

Every $j\in J\cup\left\{  1\right\}  $ satisfies
\begin{equation}
j\leq n \label{pf.prop.S(infty).inversions.a.1}%
\end{equation}
(since $n$ is the greatest element of $J\cup\left\{  1\right\}  $). Now, let
$i$ be an integer such that $i>n$. Then, $i>n>0$, so that $i$ is a positive
integer. If we had $i\in J$, then we would have $i\in J\subseteq J\cup\left\{
1\right\}  $ and thus $i\leq n$ (by (\ref{pf.prop.S(infty).inversions.a.1}),
applied to $j=i$), which would contradict $i>n$. Hence, we cannot have $i\in
J$. We thus have $i\notin J$. Since $i\in\left\{  1,2,3,\ldots\right\}  $,
this shows that $i\in\left\{  1,2,3,\ldots\right\}  \setminus J$. Thus,
$\sigma\left(  i\right)  =i$ (by (\ref{pf.prop.S(infty).inversions.a.J})).

Let us now forget that we fixed $i$. We thus have shown that every integer
$i>n$ satisfies $\sigma\left(  i\right)  =i$. Hence, there exists some
$N\in\left\{  1,2,3,\ldots\right\}  $ such that every integer $i>N$ satisfies
$\sigma\left(  i\right)  =i$ (namely, $N=n$). Proposition
\ref{prop.S(infty).inversions} \textbf{(a)} is thus proven.

\textbf{(b)} Proposition \ref{prop.S(infty).inversions} \textbf{(a)} shows
that there exists some $N\in\left\{  1,2,3,\ldots\right\}  $ such that%
\begin{equation}
\text{every integer }i>N\text{ satisfies }\sigma\left(  i\right)  =i.
\label{pf.prop.S(infty).inversions.b.N}%
\end{equation}
Consider such an $N$. We shall now show that
\[
\text{every inversion of }\sigma\text{ is an element of }\left\{
1,2,\ldots,N\right\}  ^{2}.
\]

In fact, let $c$ be an inversion of $\sigma$. We will show that $c$ is an
element of $\left\{  1,2,\ldots,N\right\}  ^{2}$.

We know that $c$ is an inversion of $\sigma$. In other words, $c$ is a pair
$\left(  i,j\right)  $ of integers satisfying $1\leq i<j$ and $\sigma\left(
i\right)  >\sigma\left(  j\right)  $ (by the definition of an
\textquotedblleft inversion of $\sigma$\textquotedblright). Consider this
$\left(  i,j\right)  $. We then have $i\leq N$\ \ \ \footnote{\textit{Proof.}
Assume the contrary. Thus, $i>N$. Hence,
(\ref{pf.prop.S(infty).inversions.b.N}) shows that $\sigma\left(  i\right)
=i$. Also, $i<j$, so that $j>i>N$. Hence,
(\ref{pf.prop.S(infty).inversions.b.N}) (applied to $j$ instead of $i$) shows
that $\sigma\left(  j\right)  =j$. Thus, $\sigma\left(  i\right)
=i<j=\sigma\left(  j\right)  $. This contradicts $\sigma\left(  i\right)
>\sigma\left(  j\right)  $. This contradiction shows that our assumption was
wrong, qed.} and $j\leq N$\ \ \ \ \footnote{\textit{Proof.} Assume the
contrary. Thus, $j>N$. Hence, (\ref{pf.prop.S(infty).inversions.b.N}) (applied
to $j$ instead of $i$) shows that $\sigma\left(  j\right)  =j$. Now,
$\sigma\left(  i\right)  >\sigma\left(  j\right)  =j>N$. Therefore,
(\ref{pf.prop.S(infty).inversions.b.N}) (applied to $\sigma\left(  i\right)  $
instead of $i$) yields $\sigma\left(  \sigma\left(  i\right)  \right)
=\sigma\left(  i\right)  $. But $\sigma$ is a permutation, and thus an
injective map. Hence, from $\sigma\left(  \sigma\left(  i\right)  \right)
=\sigma\left(  i\right)  $, we obtain $\sigma\left(  i\right)  =i$. Thus,
$\sigma\left(  i\right)  =i<j=\sigma\left(  j\right)  $. This contradicts
$\sigma\left(  i\right)  >\sigma\left(  j\right)  $. This contradiction shows
that our assumption was wrong, qed.}. From $1\leq i\leq N$, we obtain
$i\in\left\{  1,2,\ldots,N\right\}  $. From $1\leq j\leq N$, we obtain
$j\in\left\{  1,2,\ldots,N\right\}  $. Combining $i\in\left\{  1,2,\ldots
,N\right\}  $ with $j\in\left\{  1,2,\ldots,N\right\}  $, we obtain $\left(
i,j\right)  \in\left\{  1,2,\ldots,N\right\}  ^{2}$. Hence, $c=\left(
i,j\right)  \in\left\{  1,2,\ldots,N\right\}  ^{2}$. In other words, $c$ is an
element of $\left\{  1,2,\ldots,N\right\}  ^{2}$.

Now, let us forget that we fixed $c$. We thus have shown that if $c$ is an
inversion of $\sigma$, then $c$ is an element of $\left\{  1,2,\ldots
,N\right\}  ^{2}$. In other words, every inversion of $\sigma$ is an element
of $\left\{  1,2,\ldots,N\right\}  ^{2}$. Thus, there are only finitely many
inversions of $\sigma$ (since there are only finitely many elements of
$\left\{  1,2,\ldots,N\right\}  ^{2}$). Proposition
\ref{prop.S(infty).inversions} \textbf{(b)} is thus proven.
\end{proof}
\end{verlong}

Actually, Proposition \ref{prop.S(infty).inversions} \textbf{(b)} has a
converse: If a permutation $\sigma\in S_{\infty}$ has only finitely many
inversions, then $\sigma$ belongs to $S_{\left(  \infty\right)  }$. This is
easy to prove; but we won't use this.

If $\sigma\in S_{\left(  \infty\right)  }$ is a permutation, then the
\textit{length} of $\sigma$ means the number of inversions of $\sigma$. This
is well-defined, because there are only finitely many inversions of $\sigma$
(by Proposition \ref{prop.S(infty).inversions} \textbf{(b)}). The length of
$\sigma$ is denoted by $\ell\left(  \sigma\right)  $; it is a nonnegative
integer. The only permutation having length $0$ is the identity permutation
$\operatorname*{id}\in S_{\infty}$.

We have the following analogue of Exercise \ref{exe.ps2.2.5}:

\begin{exercise}
\label{exe.ps2.2.5'}\textbf{(a)} Show that every permutation $\sigma\in
S_{\left(  \infty\right)  }$ and every $k\in\left\{  1,2,3,\ldots\right\}  $
satisfy%
\begin{equation}
\ell\left(  \sigma\circ s_{k}\right)  =%
\begin{cases}
\ell\left(  \sigma\right)  +1, & \text{if }\sigma\left(  k\right)
<\sigma\left(  k+1\right)  ;\\
\ell\left(  \sigma\right)  -1, & \text{if }\sigma\left(  k\right)
>\sigma\left(  k+1\right)
\end{cases}
\label{eq.exe.2.5'.a.1}%
\end{equation}
and%
\begin{equation}
\ell\left(  s_{k}\circ\sigma\right)  =%
\begin{cases}
\ell\left(  \sigma\right)  +1, & \text{if }\sigma^{-1}\left(  k\right)
<\sigma^{-1}\left(  k+1\right)  ;\\
\ell\left(  \sigma\right)  -1, & \text{if }\sigma^{-1}\left(  k\right)
>\sigma^{-1}\left(  k+1\right)
\end{cases}
. \label{eq.exe.2.5'.a.2}%
\end{equation}

\textbf{(b)} Show that any two permutations $\sigma$ and $\tau$ in $S_{\left(
\infty\right)  }$ satisfy $\ell\left(  \sigma\circ\tau\right)  \equiv
\ell\left(  \sigma\right)  +\ell\left(  \tau\right)  \operatorname{mod}2$.

\textbf{(c)} Show that any two permutations $\sigma$ and $\tau$ in $S_{\left(
\infty\right)  }$ satisfy $\ell\left(  \sigma\circ\tau\right)  \leq\ell\left(
\sigma\right)  +\ell\left(  \tau\right)  $.

\textbf{(d)} If $\sigma\in S_{\left(  \infty\right)  }$ is a permutation
satisfying $\sigma\left(  1\right)  \leq\sigma\left(  2\right)  \leq
\sigma\left(  3\right)  \leq\cdots$, then show that $\sigma=\operatorname*{id}%
$.

\textbf{(e)} Let $\sigma\in S_{\left(  \infty\right)  }$. Show that $\sigma$
can be written as a composition of $\ell\left(  \sigma\right)  $ permutations
of the form $s_{k}$ (with $k\in\left\{  1,2,3,\ldots\right\}  $).

\textbf{(f)} Let $\sigma\in S_{\left(  \infty\right)  }$. Then, show that
$\ell\left(  \sigma\right)  =\ell\left(  \sigma^{-1}\right)  $.

\textbf{(g)} Let $\sigma\in S_{\left(  \infty\right)  }$. Show that
$\ell\left(  \sigma\right)  $ is the smallest $N\in\mathbb{N}$ such that
$\sigma$ can be written as a composition of $N$ permutations of the form
$s_{k}$ (with $k\in\left\{  1,2,3,\ldots\right\}  $).
\end{exercise}

We also have an analogue of Exercise \ref{exe.ps2.2.6}:

\begin{exercise}
\label{exe.ps2.2.6'}Let $\sigma\in S_{\left(  \infty\right)  }$. In Exercise
\ref{exe.ps2.2.4'} \textbf{(b)}, we have seen that $\sigma$ can be written as
a composition of several permutations of the form $s_{k}$ (with $k\in\left\{
1,2,3,\ldots\right\}  $). Usually there will be several ways to do so (for
instance, $\operatorname*{id}=s_{1}\circ s_{1}=s_{2}\circ s_{2}=s_{3}\circ
s_{3}=\cdots$). Show that, whichever of these ways we take, the number of
permutations composed will be congruent to $\ell\left(  \sigma\right)  $
modulo $2$.
\end{exercise}

Having defined the length of a permutation $\sigma\in S_{\left(
\infty\right)  }$, we can now define the sign of such a permutation. Again, we
mimic the definition of the sign of a $\sigma\in S_{n}$:

\begin{definition}
\label{def.perm.sign'}We define the \textit{sign} of a permutation $\sigma\in
S_{\left(  \infty\right)  }$ as the integer $\left(  -1\right)  ^{\ell\left(
\sigma\right)  }$. We denote this sign by $\left(  -1\right)  ^{\sigma}$ or
$\operatorname*{sign}\sigma$ or $\operatorname*{sgn}\sigma$. We say that a
permutation $\sigma$ is \textit{even} if its sign is $1$ (that is, if
$\ell\left(  \sigma\right)  $ is even), and \textit{odd} if its sign is $-1$
(that is, if $\ell\left(  \sigma\right)  $ is odd).
\end{definition}

Signs of permutations have the following properties:

\begin{proposition}
\label{prop.perm.signs'.basics}\textbf{(a)} The sign of the identity
permutation $\operatorname*{id}\in S_{\left(  \infty\right)  }$ is $\left(
-1\right)  ^{\operatorname*{id}}=1$. In other words, $\operatorname*{id}\in
S_{\left(  \infty\right)  }$ is even.

\textbf{(b)} For every $k\in\left\{  1,2,3,\ldots\right\}  $, the sign of the
permutation $s_{k}\in S_{\left(  \infty\right)  }$ is $\left(  -1\right)
^{s_{k}}=-1$.

\textbf{(c)} If $\sigma$ and $\tau$ are two permutations in $S_{\left(
\infty\right)  }$, then $\left(  -1\right)  ^{\sigma\circ\tau}=\left(
-1\right)  ^{\sigma}\cdot\left(  -1\right)  ^{\tau}$.

\textbf{(d)} If $\sigma\in S_{\left(  \infty\right)  }$, then $\left(
-1\right)  ^{\sigma^{-1}}=\left(  -1\right)  ^{\sigma}$.
\end{proposition}

The proof of Proposition \ref{prop.perm.signs'.basics} is analogous to the
proof of Proposition \ref{prop.perm.signs.basics}.

\begin{remark}
We have defined the sign of a permutation $\sigma\in S_{\left(  \infty\right)
}$. No such notion exists for permutations $\sigma\in S_{\infty}$. In fact,
one can show that if an element $\lambda_{\sigma}$ of $\left\{  1,-1\right\}
$ is chosen for each $\sigma\in S_{\infty}$ in such a way that every two
permutations $\sigma,\tau\in S_{\infty}$ satisfy $\lambda_{\sigma\circ\tau
}=\lambda_{\sigma}\cdot\lambda_{\tau}$, then all of the $\lambda_{\sigma}$ are
$1$. (Indeed, this follows from a result of Oystein Ore; see \url{http://mathoverflow.net/questions/54371}\ .)
\end{remark}

\begin{remark}
\label{rmk.perm.inf.lazy-ext}For every $n\in\mathbb{N}$ and every $\sigma\in
S_{n}$, we can define a permutation $\sigma_{\left(  \infty\right)  }\in
S_{\left(  \infty\right)  }$ by setting%
\[
\sigma_{\left(  \infty\right)  }\left(  i\right)  =%
\begin{cases}
\sigma\left(  i\right)  , & \text{if }i\leq n;\\
i, & \text{if }i>n
\end{cases}
\ \ \ \ \ \ \ \ \ \ \text{for all }i\in\left\{  1,2,3,\ldots\right\}  .
\]
Essentially, $\sigma_{\left(  \infty\right)  }$ is the permutation $\sigma$
extended to the set of all positive integers in the laziest possible way: It
just sends each $i>n$ to itself.

For every $n\in\mathbb{N}$, there is an injective map $\mathbf{i}_{n}%
:S_{n}\rightarrow S_{\left(  \infty\right)  }$ defined as follows:%
\[
\mathbf{i}_{n}\left(  \sigma\right)  =\sigma_{\left(  \infty\right)
}\ \ \ \ \ \ \ \ \ \ \text{for every }\sigma\in S_{n}.
\]
This map $\mathbf{i}_{n}$ is an example of what algebraists call a
\textit{group homomorphism}: It satisfies%
\begin{align*}
\mathbf{i}_{n}\left(  \operatorname*{id}\right)   &  =\operatorname*{id};\\
\mathbf{i}_{n}\left(  \sigma\circ\tau\right)   &  =\mathbf{i}_{n}\left(
\sigma\right)  \circ\mathbf{i}_{n}\left(  \tau\right)
\ \ \ \ \ \ \ \ \ \ \text{for all }\sigma,\tau\in S_{n};\\
\mathbf{i}_{n}\left(  \sigma^{-1}\right)   &  =\left(  \mathbf{i}_{n}\left(
\sigma\right)  \right)  ^{-1}\ \ \ \ \ \ \ \ \ \ \text{for all }\sigma\in
S_{n}.
\end{align*}
(This is all easy to check.) Thus, we can consider the image $\mathbf{i}%
_{n}\left(  S_{n}\right)  $ of $S_{n}$ under this map as a \textquotedblleft
copy\textquotedblright\ of $S_{n}$ which is \textquotedblleft just as good as
the original\textquotedblright\ (i.e., composition in this copy behaves in the
same way as composition in the original). It is easy to characterize this copy
inside $S_{\left(  \infty\right)  }$: Namely,%
\[
\mathbf{i}_{n}\left(  S_{n}\right)  =\left\{  \sigma\in S_{\left(
\infty\right)  }\ \mid\ \sigma\left(  i\right)  =i\text{ for all }i>n\right\}
.
\]

Using Proposition \ref{prop.S(infty).inversions} \textbf{(a)}, it is easy to
check that $S_{\left(  \infty\right)  }=\bigcup_{n\in\mathbb{N}}\mathbf{i}%
_{n}\left(  S_{n}\right)  =\mathbf{i}_{0}\left(  S_{0}\right)  \cup
\mathbf{i}_{1}\left(  S_{1}\right)  \cup\mathbf{i}_{2}\left(  S_{2}\right)
\cup\cdots$. Therefore, many properties of $S_{\left(  \infty\right)  }$ can
be derived from analogous properties of $S_{n}$ for finite $n$. For example,
using this tactic, we could easily derive Exercise \ref{exe.ps2.2.5'} from
Exercise \ref{exe.ps2.2.5}, and derive Exercise \ref{exe.ps2.2.6'} from
Exercise \ref{exe.ps2.2.6}. (However, we can just as well solve Exercises
\ref{exe.ps2.2.5'} and \ref{exe.ps2.2.6'} by copying the solutions of
Exercises \ref{exe.ps2.2.5} and \ref{exe.ps2.2.6} and making the necessary
changes; this is how I solve these exercises further below.)
\end{remark}

\subsection{\label{sect.perm.more-lengths}More on lengths of permutations}

Let us summarize some of what we have learnt about permutations. We have
defined the length $\ell\left(  \sigma\right)  $ and the inversions of a
permutation $\sigma\in S_{n}$, where $n$ is a nonnegative integer. We recall
the basic properties of these objects:

\begin{itemize}
\item For each $k\in\left\{  1,2,\ldots,n-1\right\}  $, we defined $s_{k}$ to
be the permutation in $S_{n}$ that swaps $k$ with $k+1$ but leaves all other
numbers unchanged. These permutations satisfy $s_{i}^{2}=\operatorname*{id}$
for every $i\in\left\{  1,2,\ldots,n-1\right\}  $ and%
\begin{equation}
s_{i}\circ s_{i+1}\circ s_{i}=s_{i+1}\circ s_{i}\circ s_{i+1}%
\ \ \ \ \ \ \ \ \ \ \text{for all }i\in\left\{  1,2,\ldots,n-2\right\}
\label{eq.perms.braid3}%
\end{equation}
(according to Exercise \ref{exe.ps2.2.4} \textbf{(a)}). Also, it is easy to
check that%
\begin{equation}
s_{i}\circ s_{j}=s_{j}\circ s_{i}\ \ \ \ \ \ \ \ \ \ \text{for all }%
i,j\in\left\{  1,2,\ldots,n-1\right\}  \text{ with }\left\vert i-j\right\vert
>1. \label{eq.perms.braid2}%
\end{equation}

\item An \textit{inversion} of a permutation $\sigma\in S_{n}$ means a pair
$\left(  i,j\right)  $ of integers satisfying $1\leq i<j\leq n$ and
$\sigma\left(  i\right)  >\sigma\left(  j\right)  $. The \textit{length}
$\ell\left(  \sigma\right)  $ of a permutation $\sigma\in S_{n}$ is the number
of inversions of $\sigma$.

\item Any two permutations $\sigma\in S_{n}$ and $\tau\in S_{n}$ satisfy%
\begin{equation}
\ell\left(  \sigma\circ\tau\right)  \equiv\ell\left(  \sigma\right)
+\ell\left(  \tau\right)  \operatorname{mod}2 \label{eq.perms.length.mod2}%
\end{equation}
(by Exercise \ref{exe.ps2.2.5} \textbf{(b)}) and%
\begin{equation}
\ell\left(  \sigma\circ\tau\right)  \leq\ell\left(  \sigma\right)
+\ell\left(  \tau\right)  \label{eq.perms.length.leq}%
\end{equation}
(by Exercise \ref{exe.ps2.2.5} \textbf{(c)}).

\item If $\sigma\in S_{n}$, then $\ell\left(  \sigma\right)  =\ell\left(
\sigma^{-1}\right)  $ (according to Exercise \ref{exe.ps2.2.5} \textbf{(f)}).

\item If $\sigma\in S_{n}$, then $\ell\left(  \sigma\right)  $ is the smallest
$N\in\mathbb{N}$ such that $\sigma$ can be written as a composition of $N$
permutations of the form $s_{k}$ (with $k\in\left\{  1,2,\ldots,n-1\right\}
$). (This follows from Exercise \ref{exe.ps2.2.5} \textbf{(g)}.)
\end{itemize}

By now, we know almost all about the $s_{k}$'s and about the lengths of
permutations that is necessary for studying determinants. (\textquotedblleft
Almost\textquotedblright\ because Exercise \ref{exe.ps4.1ab} below will also
be useful.) I shall now sketch some more advanced properties of these things,
partly as a curiosity, partly to further your intuition; none of these
properties shall be used further below. The rest of Section
\ref{sect.perm.more-lengths} shall rely on some notions we have not introduced
in these notes; in particular, we will use the concepts of undirected graphs
(\cite[Chapter 12]{LeLeMe16}), directed graphs (\cite[Chapter 10]{LeLeMe16})
and (briefly) polytopes (see, e.g., \cite[Chapter 10]{AigZie}).

First, here is a way to visualize lengths of permutations using graph theory:

Fix $n\in\mathbb{N}$. We define the $n$\textit{-th right Bruhat graph} to be
the (undirected) graph whose vertices are the permutations $\sigma\in S_{n}$,
and whose edges are defined as follows: Two vertices $\sigma\in S_{n}$ and
$\tau\in S_{n}$ are adjacent if and only if there exists a $k\in\left\{
1,2,\ldots,n-1\right\}  $ such that $\sigma=\tau\circ s_{k}$. (This condition
is clearly symmetric in $\sigma$ and $\tau$: If $\sigma=\tau\circ s_{k}$, then
$\tau=\sigma\circ s_{k}$.) For instance, the $3$-rd right Bruhat graph looks
as follows:%
\[%
\xymatrix{
& 321 \ar@{-}[dl] \ar@{-}[dr] & \\
312 \ar@{-}[d] & & 231 \ar@{-}[d] \\
132 \ar@{-}[dr] & & 213 \ar@{-}[dl] \\
& 123 &
}%
,
\]
where we are writing permutations in one-line notation (and omitting
parentheses and commas). The $4$-th right Bruhat graph can be seen
\href{https://en.wikipedia.org/wiki/File:Symmetric_group_4;_Cayley_graph_1,2,6_%283D%29.svg}{on
Wikipedia}.\footnote{Don't omit the word \textquotedblleft
right\textquotedblright\ in \textquotedblleft right Bruhat
graph\textquotedblright; else it means a different graph with more edges.}

These graphs have lots of properties. There is a canonical way to direct their
edges: The edge between $\sigma$ and $\tau$ is directed towards the vertex
with the larger length. (The lengths of $\sigma$ and $\tau$ always differ by
$1$ if there is an edge between $\sigma$ and $\tau$.) This way, the $n$-th
right Bruhat graph is an acyclic directed graph. It therefore defines a
partially ordered set, called the \textit{right permutohedron order}%
\footnote{also known as the \textit{right weak order} or \textit{right weak
Bruhat order} (but, again, do not omit the words \textquotedblleft
right\textquotedblright\ and \textquotedblleft weak\textquotedblright)} on
$S_{n}$, whose elements are the permutations $\sigma\in S_{n}$ and whose order
relation is defined as follows: We have $\sigma\leq\tau$ if and only if there
is a directed path from $\sigma$ to $\tau$ in the directed $n$-th right Bruhat
graph. If you know
\href{https://en.wikipedia.org/wiki/Lattice_%28order%29}{the (combinatorial)
notion of a lattice}, you might notice that this right permutohedron order
\href{http://mathoverflow.net/questions/158945/why-is-the-right-permutohedron-order-aka-weak-order-on-s-n-a-lattice}{is
a lattice}.

The word \textquotedblleft permutohedron\textquotedblright\ in
\textquotedblleft permutohedron order\textquotedblright\ hints at what might
be its least expected property: The $n$-th Bruhat graph can be viewed as the
graph formed by the vertices and the edges of a certain polytope in
$n$-dimensional space $\mathbb{R}^{n}$. This polytope -- called the
$n$\textit{-th
\href{https://en.wikipedia.org/wiki/Permutohedron}{\textit{permutohedron}}%
}\footnote{Some spell it \textquotedblleft\textit{permutahedron}%
\textquotedblright\ instead. The word is of relatively recent origin (1969).}
-- is the convex hull of the points $\left(  \sigma\left(  1\right)
,\sigma\left(  2\right)  ,\ldots,\sigma\left(  n\right)  \right)  $ for
$\sigma\in S_{n}$. These points are its vertices; however, its vertex $\left(
\sigma\left(  1\right)  ,\sigma\left(  2\right)  ,\ldots,\sigma\left(
n\right)  \right)  $ corresponds to the vertex $\sigma^{-1}$ (not $\sigma$) of
the $n$-th Bruhat graph. Feel free to roam
\href{https://en.wikipedia.org/wiki/Permutohedron}{its Wikipedia page} for
other (combinatorial and geometric) curiosities.

The notion of a length fits perfectly into this whole picture. For instance,
the length $\ell\left(  \sigma\right)  $ of a permutation $\sigma$ is the
smallest length of a path from $\operatorname*{id}\in S_{n}$ to $\sigma$ on
the $n$-th right Bruhat graph (and this holds no matter whether the graph is
considered to be directed or not). For the undirected Bruhat graphs, we have
something more general:

\begin{exercise}
\label{exe.ps4.0}Let $n\in\mathbb{N}$. Let $\sigma\in S_{n}$ and $\tau\in
S_{n}$. Show that $\ell\left(  \sigma^{-1}\circ\tau\right)  $ is the smallest
length of a path between $\sigma$ and $\tau$ on the (undirected) $n$-th right
Bruhat graph.
\end{exercise}

(Recall that the length of a path in a graph is defined as the number of edges
along this path.)

How many permutations in $S_{n}$ have a given length? The number is not easy
to compute directly; however, its generating function is nice. (See
\cite[Chapter 16]{LeLeMe16} for the notion of generating functions.) Namely,%
\[
\sum_{w\in S_{n}}q^{\ell\left(  w\right)  }=\left(  1+q\right)  \left(
1+q+q^{2}\right)  \cdots\left(  1+q+q^{2}+\cdots+q^{n-1}\right)
\]
(where $q$ is an indeterminate). This equality (with $q$ renamed as $x$) is
Corollary~\ref{cor.perm.lehmer.lensum} (which is proven below, in the solution
to Exercise~\ref{exe.perm.lehmer.prove}). Another proof can be found in
\cite[Corollary 1.3.13]{Stanley-EC1} (but notice that Stanley denotes $S_{n}$
by $\mathfrak{S}_{n}$, and $\ell\left(  w\right)  $ by $\operatorname*{inv}%
\left(  w\right)  $).

\begin{remark}
Much more can be said. Let me briefly mention (without proof) two other
related results.

We can ask ourselves in what different ways a permutation can be written as a
composition of $N$ permutations of the form $s_{k}$. For instance, the
permutation $w_{0}\in S_{3}$ which sends $1$, $2$ and $3$ to $3$, $2$ and $1$,
respectively (that is, $w_{0}=\left(  3,2,1\right)  $ in one-line notation)
can be written as a product of three $s_{k}$'s in the two forms%
\begin{equation}
w_{0}=s_{1}\circ s_{2}\circ s_{1},\ \ \ \ \ \ \ \ \ \ w_{0}=s_{2}\circ
s_{1}\circ s_{2}, \label{eq.tits.example}%
\end{equation}
but can also be written as a product of five $s_{k}$'s (e.g., as $w_{0}%
=s_{1}\circ s_{2}\circ s_{1}\circ s_{2}\circ s_{2}$) or seven $s_{k}$'s or
nine $s_{k}$'s, etc. Are the different representations of $w_{0}$ related?

Clearly, the two representations in (\ref{eq.tits.example}) are connected to
each other by the equality $s_{1}\circ s_{2}\circ s_{1}=s_{2}\circ s_{1}\circ
s_{2}$, which is a particular case of (\ref{eq.perms.braid3}). Also, the
representation $w_{0}=s_{1}\circ s_{2}\circ s_{1}\circ s_{2}\circ s_{2}$
reduces to $w_{0}=s_{1}\circ s_{2}\circ s_{1}$ by \textquotedblleft
cancelling\textquotedblright\ the two adjacent $s_{2}$'s at the end (recall
that $s_{i}\circ s_{i}=s_{i}^{2}=\operatorname*{id}$ for every $i$).

Interestingly, this generalizes. Let $n\in\mathbb{N}$ and $\sigma\in S_{n}$. A
\textit{reduced expression} for $\sigma$ will mean a representation of
$\sigma$ as a composition of $\ell\left(  \sigma\right)  $ permutations of the
form $s_{k}$. (As we know, less than $\ell\left(  \sigma\right)  $ such
permutations do not suffice; thus the name \textquotedblleft
reduced\textquotedblright.) Then, (one of the many versions of)
\textit{Matsumoto's theorem} states that any two reduced expressions of
$\sigma$ can be obtained from one another by a rewriting process, each step of
which is either an application of (\ref{eq.perms.braid3}) (i.e., you pick an
\textquotedblleft$s_{i}\circ s_{i+1}\circ s_{i}$\textquotedblright\ in the
expression and replace it by \textquotedblleft$s_{i+1}\circ s_{i}\circ
s_{i+1}$\textquotedblright, or vice versa) or an application of
(\ref{eq.perms.braid2}) (i.e., you pick an \textquotedblleft$s_{i}\circ s_{j}%
$\textquotedblright\ with $\left\vert i-j\right\vert >1$ and replace it by
\textquotedblleft$s_{j}\circ s_{i}$\textquotedblright, or vice versa). For
instance, for $n=4$ and $\sigma=\left(  4,3,1,2,5\right)  $ (in one-line
notation), the two reduced expressions $\sigma=s_{1}\circ s_{2}\circ
s_{3}\circ s_{1}\circ s_{2}$ and $\sigma=s_{2}\circ s_{3}\circ s_{1}\circ
s_{2}\circ s_{3}$ can be obtained from one another by the following rewriting
process:%
\begin{align*}
s_{1}\circ s_{2}\circ\underbrace{s_{3}\circ s_{1}}_{\substack{=s_{1}\circ
s_{3}\\\text{(by (\ref{eq.perms.braid2}))}}}\circ s_{2}  &  =\underbrace{s_{1}%
\circ s_{2}\circ s_{1}}_{\substack{=s_{2}\circ s_{1}\circ s_{2}\\\text{(by
(\ref{eq.perms.braid3}))}}}\circ s_{3}\circ s_{2}=s_{2}\circ s_{1}%
\circ\underbrace{s_{2}\circ s_{3}\circ s_{2}}_{\substack{=s_{3}\circ
s_{2}\circ s_{3}\\\text{(by (\ref{eq.perms.braid3}))}}}\\
&  =s_{2}\circ\underbrace{s_{1}\circ s_{3}}_{\substack{=s_{3}\circ
s_{1}\\\text{(by (\ref{eq.perms.braid2}))}}}\circ s_{2}\circ s_{3}=s_{2}\circ
s_{3}\circ s_{1}\circ s_{2}\circ s_{3}.
\end{align*}
See, e.g., Williamson's thesis \cite[Corollary 1.2.3]{William03} or Knutson's
notes \cite[\S 2.3]{Knutson} for a proof of this fact. (Knutson, instead of
saying that \textquotedblleft$\sigma=s_{k_{1}}\circ s_{k_{2}}\circ\cdots\circ
s_{k_{p}}$ is a reduced expression for $\sigma$\textquotedblright, says that
\textquotedblleft$k_{1}k_{2}\cdots k_{p}$ is a reduced word for $\sigma
$\textquotedblright.)

Something subtler holds for \textquotedblleft non-reduced\textquotedblright%
\ expressions. Namely, if we have a representation of $\sigma$ as a
composition of some number of permutations of the form $s_{k}$ (not
necessarily $\ell\left(  \sigma\right)  $ of them), then we can transform it
into a reduced expression by a rewriting process which consists of
applications of (\ref{eq.perms.braid3}) and (\ref{eq.perms.braid2}) as before
and also of cancellation steps (i.e., picking an \textquotedblleft$s_{i}\circ
s_{i}$\textquotedblright\ in the expression and removing it). This follows
from \cite[Chapter SYM, Proposition (2.6)]{LLPT95}\footnotemark, and can also
easily be derived from \cite[Corollary 1.2.3 and Corollary 1.1.6]{William03}.

This all is stated and proven in greater generality in good books on Coxeter
groups, such as \cite{BjoBre05}. We won't need these results in the following,
but they are an example of what one can see if one looks at permutations closely.
\end{remark}

\footnotetext{What the authors of \cite{LLPT95} call a \textquotedblleft
presentation\textquotedblright\ of a permutation $\sigma\in S_{n}$ is a finite
list $\left(  s_{k_{1}},s_{k_{2}},\ldots,s_{k_{p}}\right)  $ of elements of
$\left\{  s_{1},s_{2},\ldots,s_{n-1}\right\}  $ satisfying $\sigma=s_{k_{1}%
}\circ s_{k_{2}}\circ\cdots\circ s_{k_{p}}$. What the authors of \cite{LLPT95}
call a \textquotedblleft minimal presentation\textquotedblright\ of $\sigma$
is what we call a reduced expression of $\sigma$.}

\subsection{More on signs of permutations}

In Section \ref{sect.sign}, we have defined the sign $\left(  -1\right)
^{\sigma}=\operatorname*{sign}\sigma=\operatorname*{sgn}\sigma$ of a
permutation $\sigma$. We recall the most important facts about it:

\begin{itemize}
\item We have $\left(  -1\right)  ^{\sigma}=\left(  -1\right)  ^{\ell\left(
\sigma\right)  }$ for every $\sigma\in S_{n}$. (This is the definition of
$\left(  -1\right)  ^{\sigma}$.) Thus, for every $\sigma\in S_{n}$, we have
$\left(  -1\right)  ^{\sigma}=\left(  -1\right)  ^{\ell\left(  \sigma\right)
}\in\left\{  1,-1\right\}  $.

\item The permutation $\sigma\in S_{n}$ is said to be \textit{even} if
$\left(  -1\right)  ^{\sigma}=1$, and is said to be \textit{odd} if $\left(
-1\right)  ^{\sigma}=-1$.

\item The sign of the identity permutation $\operatorname*{id}\in S_{n}$ is
$\left(  -1\right)  ^{\operatorname*{id}}=1$.

\item For every $k\in\left\{  1,2,\ldots,n-1\right\}  $, the sign of the
permutation $s_{k}\in S_{n}$ is $\left(  -1\right)  ^{s_{k}}=-1$.

\item If $\sigma$ and $\tau$ are two permutations in $S_{n}$, then
\begin{equation}
\left(  -1\right)  ^{\sigma\circ\tau}=\left(  -1\right)  ^{\sigma}\cdot\left(
-1\right)  ^{\tau}. \label{eq.sign.prod}%
\end{equation}

\item If $\sigma\in S_{n}$, then%
\begin{equation}
\left(  -1\right)  ^{\sigma^{-1}}=\left(  -1\right)  ^{\sigma}.
\label{eq.sign.inverse}%
\end{equation}

\end{itemize}

A simple consequence of the above facts is the following proposition:

\begin{proposition}
\label{prop.sign.prod-of-many}Let $n\in\mathbb{N}$ and $k\in\mathbb{N}$. Let
$\sigma_{1},\sigma_{2},\ldots,\sigma_{k}$ be $k$ permutations in $S_{n}$.
Then,%
\begin{equation}
\left(  -1\right)  ^{\sigma_{1}\circ\sigma_{2}\circ\cdots\circ\sigma_{k}%
}=\left(  -1\right)  ^{\sigma_{1}}\cdot\left(  -1\right)  ^{\sigma_{2}}%
\cdot\cdots\cdot\left(  -1\right)  ^{\sigma_{k}}. \label{eq.sign.prod-of-many}%
\end{equation}

\end{proposition}

\begin{vershort}
\begin{proof}
[Proof of Proposition \ref{prop.sign.prod-of-many}.]Straightforward induction
over $k$. The induction base (i.e., the case when $k=0$) follows from the fact
that $\left(  -1\right)  ^{\operatorname*{id}}=1$ (since the composition of
$0$ permutations is $\operatorname*{id}$). The induction step is easily done
using (\ref{eq.sign.prod}).
\end{proof}
\end{vershort}

\begin{verlong}
\begin{proof}
[Proof of Proposition \ref{prop.sign.prod-of-many}.]We will prove
(\ref{eq.sign.prod-of-many}) by induction over $k$:

\textit{Induction base:} Assume that $k=0$. Then, $\sigma_{1}\circ\sigma
_{2}\circ\cdots\circ\sigma_{k}=\sigma_{1}\circ\sigma_{2}\circ\cdots\circ
\sigma_{0}=\operatorname*{id}$ (since the composition of $0$ permutations is
$\operatorname*{id}$). Hence, $\left(  -1\right)  ^{\sigma_{1}\circ\sigma
_{2}\circ\cdots\circ\sigma_{k}}=\left(  -1\right)  ^{\operatorname*{id}}=1$.
On the other hand, from $k=0$, we obtain $\left(  -1\right)  ^{\sigma_{1}%
}\cdot\left(  -1\right)  ^{\sigma_{2}}\cdot\cdots\cdot\left(  -1\right)
^{\sigma_{k}}=\left(  -1\right)  ^{\sigma_{1}}\cdot\left(  -1\right)
^{\sigma_{2}}\cdot\cdots\cdot\left(  -1\right)  ^{\sigma_{0}}=1$ (since the
product of $0$ integers is $1$). Compared with $\left(  -1\right)
^{\sigma_{1}\circ\sigma_{2}\circ\cdots\circ\sigma_{k}}=1$, this yields
$\left(  -1\right)  ^{\sigma_{1}\circ\sigma_{2}\circ\cdots\circ\sigma_{k}%
}=\left(  -1\right)  ^{\sigma_{1}}\cdot\left(  -1\right)  ^{\sigma_{2}}%
\cdot\cdots\cdot\left(  -1\right)  ^{\sigma_{k}}$. Thus,
(\ref{eq.sign.prod-of-many}) is proven for $k=0$. The induction base is thus complete.

\textit{Induction step:} Let $K\in\mathbb{N}$. Assume that
(\ref{eq.sign.prod-of-many}) is proven for $k=K$. We need to prove
(\ref{eq.sign.prod-of-many}) for $k=K+1$.

Let $\sigma_{1},\sigma_{2},\ldots,\sigma_{K+1}$ be $K+1$ permutations in
$S_{n}$. We have assumed that (\ref{eq.sign.prod-of-many}) is proven for
$k=K$. Thus, we can apply (\ref{eq.sign.prod-of-many}) to $k=K$. We thus
conclude $\left(  -1\right)  ^{\sigma_{1}\circ\sigma_{2}\circ\cdots\circ
\sigma_{K}}=\left(  -1\right)  ^{\sigma_{1}}\cdot\left(  -1\right)
^{\sigma_{2}}\cdot\cdots\cdot\left(  -1\right)  ^{\sigma_{K}}$.

But $\sigma_{1}\circ\sigma_{2}\circ\cdots\circ\sigma_{K+1}=\left(  \sigma
_{1}\circ\sigma_{2}\circ\cdots\circ\sigma_{K}\right)  \circ\sigma_{K+1}$, so
that%
\begin{align*}
\left(  -1\right)  ^{\sigma_{1}\circ\sigma_{2}\circ\cdots\circ\sigma_{K+1}}
&  =\left(  -1\right)  ^{\left(  \sigma_{1}\circ\sigma_{2}\circ\cdots
\circ\sigma_{K}\right)  \circ\sigma_{K+1}}=\underbrace{\left(  -1\right)
^{\sigma_{1}\circ\sigma_{2}\circ\cdots\circ\sigma_{K}}}_{=\left(  -1\right)
^{\sigma_{1}}\cdot\left(  -1\right)  ^{\sigma_{2}}\cdot\cdots\cdot\left(
-1\right)  ^{\sigma_{K}}}\cdot\left(  -1\right)  ^{\sigma_{K+1}}\\
&  \ \ \ \ \ \ \ \ \ \ \left(  \text{by (\ref{eq.sign.prod}), applied to
}\sigma=\sigma_{1}\circ\sigma_{2}\circ\cdots\circ\sigma_{K}\text{ and }%
\tau=\sigma_{K+1}\right) \\
&  =\left(  \left(  -1\right)  ^{\sigma_{1}}\cdot\left(  -1\right)
^{\sigma_{2}}\cdot\cdots\cdot\left(  -1\right)  ^{\sigma_{K}}\right)
\cdot\left(  -1\right)  ^{\sigma_{K+1}}\\
&  =\left(  -1\right)  ^{\sigma_{1}}\cdot\left(  -1\right)  ^{\sigma_{2}}%
\cdot\cdots\cdot\left(  -1\right)  ^{\sigma_{K+1}}.
\end{align*}

Let us now forget that we fixed $\sigma_{1},\sigma_{2},\ldots,\sigma_{K+1}$.
We thus have shown that if $\sigma_{1},\sigma_{2},\ldots,\sigma_{K+1}$ are
$K+1$ permutations in $S_{n}$, then $\left(  -1\right)  ^{\sigma_{1}%
\circ\sigma_{2}\circ\cdots\circ\sigma_{K+1}}=\left(  -1\right)  ^{\sigma_{1}%
}\cdot\left(  -1\right)  ^{\sigma_{2}}\cdot\cdots\cdot\left(  -1\right)
^{\sigma_{K+1}}$. In other words, we have proven (\ref{eq.sign.prod-of-many})
for $k=K+1$. This completes the induction step, and so the induction proof of
(\ref{eq.sign.prod-of-many}) is complete. In other words, Proposition
\ref{prop.sign.prod-of-many} is proven.
\end{proof}
\end{verlong}

Let us introduce another notation:

\begin{definition}
\label{def.transpos}Let $n\in\mathbb{N}$. Let $i$ and $j$ be two distinct
elements of $\left\{  1,2,\ldots,n\right\}  $. We let $t_{i,j}$ be the
permutation in $S_{n}$ which swaps $i$ with $j$ while leaving all other
elements of $\left\{  1,2,\ldots,n\right\}  $ unchanged. Such a permutation is
called a \textit{transposition} (and is often denoted by $\left(  i,j\right)
$ in literature; but we prefer not to do so, since it is too similar to
one-line notation).
\end{definition}

Notice that the permutations $s_{1},s_{2},\ldots,s_{n-1}$ are transpositions
(namely, $s_{i}=t_{i,i+1}$ for every $i\in\left\{  1,2,\ldots,n-1\right\}  $),
but they are not the only transpositions (when $n\geq3$).

For the next exercise, we need one further definition, which extends
Definition \ref{def.transpos}:

\begin{definition}
\label{def.transpos.ii}Let $n\in\mathbb{N}$. Let $i$ and $j$ be two elements
of $\left\{  1,2,\ldots,n\right\}  $. We define a permutation $t_{i,j}\in
S_{n}$ as follows:

\begin{itemize}
\item If $i\neq j$, then the permutation $t_{i,j}$ has already been defined in
Definition \ref{def.transpos}.

\item If $i=j$, then we define the permutation $t_{i,j}$ to be the identity
$\operatorname*{id}\in S_{n}$.
\end{itemize}
\end{definition}

\begin{exercise}
\label{exe.transpos.code}Whenever $m$ is an integer, we shall use the notation
$\left[  m\right]  $ for the set $\left\{  1,2,\ldots,m\right\}  $.

Let $n\in\mathbb{N}$. Let $\sigma\in S_{n}$. Prove that there is a unique
$n$-tuple $\left(  i_{1},i_{2},\ldots,i_{n}\right)  \in\left[  1\right]
\times\left[  2\right]  \times\cdots\times\left[  n\right]  $ such that%
\[
\sigma=t_{1,i_{1}}\circ t_{2,i_{2}}\circ\cdots\circ t_{n,i_{n}}.
\]

\end{exercise}

\begin{example}
For this example, set $n=4$, and let $\sigma\in S_{4}$ be the permutation that
sends $1,2,3,4$ to $3,1,4,2$. Then, Exercise \ref{exe.transpos.code} claims
that there is a unique $4$-tuple $\left(  i_{1},i_{2},i_{3},i_{4}\right)
\in\left[  1\right]  \times\left[  2\right]  \times\left[  3\right]
\times\left[  4\right]  $ such that $\sigma=t_{1,i_{1}}\circ t_{2,i_{2}}\circ
t_{3,i_{3}}\circ t_{4,i_{4}}$.

This $4$-tuple can easily be found: it is $\left(  1,1,1,3\right)  $. In fact,
we have $\sigma=t_{1,1}\circ t_{2,1}\circ t_{3,1}\circ t_{4,3}$.
\end{example}

\begin{exercise}
\label{exe.ps4.1ab}Let $n\in\mathbb{N}$. Let $i$ and $j$ be two distinct
elements of $\left\{  1,2,\ldots,n\right\}  $.

\textbf{(a)} Find $\ell\left(  t_{i,j}\right)  $.

\textbf{(b)} Show that $\left(  -1\right)  ^{t_{i,j}}=-1$.
\end{exercise}

\begin{exercise}
\label{exe.ps4.1c}Let $n\in\mathbb{N}$. Let $w_{0}$ denote the permutation in
$S_{n}$ which sends each $k\in\left\{  1,2,\ldots,n\right\}  $ to $n+1-k$.
Compute $\ell\left(  w_{0}\right)  $ and $\left(  -1\right)  ^{w_{0}}$.
\end{exercise}

\begin{exercise}
\label{exe.ps4.2}Let $X$ be a finite set. We want to define the sign of any
permutation of $X$. (We have sketched this definition before (see
(\ref{eq.ps2.S(X).sign.teaser})), but now we shall do it in detail.)

Fix a bijection $\phi:X\rightarrow\left\{  1,2,\ldots,n\right\}  $ for some
$n\in\mathbb{N}$. (Such a bijection always exists. Indeed, constructing such a
bijection is tantamount to writing down a list of all elements of $X$, with no
duplicates.) For every permutation $\sigma$ of $X$, set%
\[
\left(  -1\right)  _{\phi}^{\sigma}=\left(  -1\right)  ^{\phi\circ\sigma
\circ\phi^{-1}}.
\]
Here, the right hand side is well-defined because $\phi\circ\sigma\circ
\phi^{-1}$ is a permutation of $\left\{  1,2,\ldots,n\right\}  $.

\textbf{(a)} Prove that $\left(  -1\right)  _{\phi}^{\sigma}$ depends only on
the permutation $\sigma$ of $X$, but not on the bijection $\phi$. (In other
words, for a given $\sigma$, any two different choices of $\phi$ will lead to
the same $\left(  -1\right)  _{\phi}^{\sigma}$.)

This allows us to define the \textit{sign} of the permutation $\sigma$ to be
$\left(  -1\right)  _{\phi}^{\sigma}$ for any choice of bijection
$\phi:X\rightarrow\left\{  1,2,\ldots,n\right\}  $. We denote this sign simply
by $\left(  -1\right)  ^{\sigma}$. (When $X=\left\{  1,2,\ldots,n\right\}  $,
then this sign is clearly the same as the sign $\left(  -1\right)  ^{\sigma}$
we defined before, because we can pick the bijection $\phi=\operatorname*{id}$.)

\textbf{(b)} Show that the permutation $\operatorname*{id}:X\rightarrow X$
satisfies $\left(  -1\right)  ^{\operatorname*{id}}=1$.

\textbf{(c)} Show that $\left(  -1\right)  ^{\sigma\circ\tau}=\left(
-1\right)  ^{\sigma}\cdot\left(  -1\right)  ^{\tau}$ for any two permutations
$\sigma$ and $\tau$ of $X$.
\end{exercise}

\begin{remark}
A sufficiently pedantic reader might have noticed that the definition of
$\left(  -1\right)  ^{\sigma}$ in Exercise \ref{exe.ps4.2} is not completely
kosher. In fact, the set $X$ may be $\left\{  1,2,\ldots,n\right\}  $ for some
$n\in\mathbb{N}$; in this case, $\sigma$ is an element of $S_{n}$, and thus
the sign $\left(  -1\right)  ^{\sigma}$ has already been defined in Definition
\ref{def.perm.sign}. Thus, in this case, we are defining the notation $\left(
-1\right)  ^{\sigma}$ a second time in Exercise \ref{exe.ps4.2}. Woe to us if
this second definition yields a different number than the first!

Fortunately, it does not. The definition of $\left(  -1\right)  ^{\sigma}$ in
Exercise \ref{exe.ps4.2} does not conflict with the original meaning of
$\left(  -1\right)  ^{\sigma}$ as defined in Definition \ref{def.perm.sign}.
Indeed, in order to prove this, we temporarily rename the number $\left(
-1\right)  ^{\sigma}$ defined in Exercise \ref{exe.ps4.2} as $\left(
-1\right)  _{\operatorname*{new}}^{\sigma}$ (in order to ensure that we don't
confuse it with the number $\left(  -1\right)  ^{\sigma}$ defined in
Definition \ref{def.perm.sign}). Now, consider the situation of Exercise
\ref{exe.ps4.2}, and assume that $X=\left\{  1,2,\ldots,n\right\}  $. We must
then prove that $\left(  -1\right)  _{\operatorname*{new}}^{\sigma}=\left(
-1\right)  ^{\sigma}$. But the definition of $\left(  -1\right)
_{\operatorname*{new}}^{\sigma}$ in Exercise \ref{exe.ps4.2} says that
$\left(  -1\right)  _{\operatorname*{new}}^{\sigma}=\left(  -1\right)  _{\phi
}^{\sigma}$, where $\phi$ is any bijection $X\rightarrow\left\{
1,2,\ldots,n\right\}  $. We can apply this to $\phi=\operatorname*{id}$
(because clearly, $\operatorname*{id}$ is a bijection $X\rightarrow\left\{
1,2,\ldots,n\right\}  $), and thus obtain $\left(  -1\right)
_{\operatorname*{new}}^{\sigma}=\left(  -1\right)  _{\operatorname*{id}%
}^{\sigma}$. But the definition of $\left(  -1\right)  _{\operatorname*{id}%
}^{\sigma}$ yields $\left(  -1\right)  _{\operatorname*{id}}^{\sigma}=\left(
-1\right)  ^{\operatorname*{id}\circ\sigma\circ\operatorname*{id}%
\nolimits^{-1}}=\left(  -1\right)  ^{\sigma}$ (since $\operatorname*{id}%
\circ\sigma\circ\underbrace{\operatorname*{id}\nolimits^{-1}}%
_{=\operatorname*{id}}=\sigma$). Thus, $\left(  -1\right)
_{\operatorname*{new}}^{\sigma}=\left(  -1\right)  _{\operatorname*{id}%
}^{\sigma}=\left(  -1\right)  ^{\sigma}$. This is precisely what we wanted to
prove. Thus, we have shown that the definition of $\left(  -1\right)
^{\sigma}$ in Exercise \ref{exe.ps4.2} does not conflict with the original
meaning of $\left(  -1\right)  ^{\sigma}$ as defined in Definition
\ref{def.perm.sign}.
\end{remark}

\begin{remark}
Let $n\in\mathbb{N}$. Recall that a \textit{transposition} in $S_{n}$ means a
permutation of the form $t_{i,j}$, where $i$ and $j$ are two distinct elements
of $\left\{  1,2,\ldots,n\right\}  $. Therefore, if $\tau$ is a transposition
in $S_{n}$, then
\begin{equation}
\left(  -1\right)  ^{\tau}=-1. \label{eq.sign.transposition}%
\end{equation}
(In fact, if $\tau$ is a transposition in $S_{n}$, then $\tau$ can be written
in the form $\tau=t_{i,j}$ for two distinct elements $i$ and $j$ of $\left\{
1,2,\ldots,n\right\}  $; and therefore, for these two elements $i$ and $j$, we
have $\left(  -1\right)  ^{\tau}=\left(  -1\right)  ^{t_{i,j}}=-1$ (according
to Exercise \ref{exe.ps4.1ab} \textbf{(b)}). This proves
(\ref{eq.sign.transposition}).)

Now, let $\sigma\in S_{n}$ be any permutation. Assume that $\sigma$ is written
in the form $\sigma=\tau_{1}\circ\tau_{2}\circ\cdots\circ\tau_{k}$ for some
transpositions $\tau_{1},\tau_{2},\ldots,\tau_{k}$ in $S_{n}$. Then,%
\begin{align}
\left(  -1\right)  ^{\sigma}  &  =\left(  -1\right)  ^{\tau_{1}\circ\tau
_{2}\circ\cdots\circ\tau_{k}}=\underbrace{\left(  -1\right)  ^{\tau_{1}}%
}_{\substack{=-1\\\text{(by (\ref{eq.sign.transposition}))}}}\cdot
\underbrace{\left(  -1\right)  ^{\tau_{2}}}_{\substack{=-1\\\text{(by
(\ref{eq.sign.transposition}))}}}\cdot\cdots\cdot\underbrace{\left(
-1\right)  ^{\tau_{k}}}_{\substack{=-1\\\text{(by (\ref{eq.sign.transposition}%
))}}}\nonumber\\
&  \ \ \ \ \ \ \ \ \ \ \left(  \text{by (\ref{eq.sign.prod-of-many}), applied
to }\sigma_{i}=\tau_{i}\right) \nonumber\\
&  =\underbrace{\left(  -1\right)  \cdot\left(  -1\right)  \cdot\cdots
\cdot\left(  -1\right)  }_{k\text{ factors}}=\left(  -1\right)  ^{k}.
\label{eq.sign.prod-of-transpositions}%
\end{align}
Since many permutations can be written as products of transpositions in a
simple way, this formula gives a useful method for computing signs.
\end{remark}

\begin{remark}
Let $n\in\mathbb{N}$. It is not hard to prove that
\begin{equation}
\left(  -1\right)  ^{\sigma}=\prod_{1\leq i<j\leq n}\dfrac{\sigma\left(
i\right)  -\sigma\left(  j\right)  }{i-j}\ \ \ \ \ \ \ \ \ \ \text{for every
}\sigma\in S_{n}. \label{eq.sign.pseudoexplicit}%
\end{equation}
(Of course, it is no easier to compute $\left(  -1\right)  ^{\sigma}$ using
this seemingly explicit formula than by counting inversions.)

We shall prove (\ref{eq.sign.pseudoexplicit}) in Exercise
\ref{exe.perm.sign.pseudoexplicit} \textbf{(c)}.
\end{remark}

\begin{remark}
The sign of a permutation is also called its \textit{signum} or its
\textit{signature}. Different authors define the sign of a permutation
$\sigma$ in different ways. Some (e.g., Hefferon in \cite[Chapter Four,
Definition 4.7]{Hefferon}, or Strickland in \cite[Definition B.4]{Strick13})
define it as we do; others (e.g., Conrad in \cite{Conrad}, or M\'{a}t\'{e} in
\cite{Mate14}, or Hoffman and Kunze in \cite[p. 152]{HoffmanKunze}) define it
using (\ref{eq.sign.prod-of-transpositions}); yet others define it using
something called the \textit{cycle decomposition} of a permutation; some even
define it using (\ref{eq.sign.pseudoexplicit}), or using a similar ratio of
two polynomials. However, it is not hard to check that all of these
definitions are equivalent. (We already know that the first two of them are equivalent.)
\end{remark}

\begin{exercise}
\label{exe.perm.sign.pseudoexplicit}Let $n\in\mathbb{N}$. Let $\sigma\in
S_{n}$.

\textbf{(a)} If $x_{1},x_{2},\ldots,x_{n}$ are $n$ elements of $\mathbb{C}$,
then prove that%
\[
\prod_{1\leq i<j\leq n}\left(  x_{\sigma\left(  i\right)  }-x_{\sigma\left(
j\right)  }\right)  =\left(  -1\right)  ^{\sigma}\cdot\prod_{1\leq i<j\leq
n}\left(  x_{i}-x_{j}\right)  .
\]

\textbf{(b)} More generally: For every $\left(  i,j\right)  \in\left\{
1,2,\ldots,n\right\}  ^{2}$, let $a_{\left(  i,j\right)  }$ be an element of
$\mathbb{C}$. Assume that%
\begin{equation}
a_{\left(  j,i\right)  }=-a_{\left(  i,j\right)  }%
\ \ \ \ \ \ \ \ \ \ \text{for every }\left(  i,j\right)  \in\left\{
1,2,\ldots,n\right\}  ^{2}. \label{eq.exe.perm.sign.pseudoexplicit.b.skew}%
\end{equation}
Prove that%
\[
\prod_{1\leq i<j\leq n}a_{\left(  \sigma\left(  i\right)  ,\sigma\left(
j\right)  \right)  }=\left(  -1\right)  ^{\sigma}\cdot\prod_{1\leq i<j\leq
n}a_{\left(  i,j\right)  }.
\]

\textbf{(c)} Prove (\ref{eq.sign.pseudoexplicit}).

\textbf{(d)} Use Exercise \ref{exe.perm.sign.pseudoexplicit} \textbf{(a)} to
give a new solution to Exercise \ref{exe.ps2.2.5} \textbf{(b)}.
\end{exercise}

The next exercise relies on the notion of ``the list of all elements of $S$ in
increasing order (with no repetitions)'', where $S$ is a finite set of
integers. This notion means exactly what it says; it was rigorously defined in
Definition \ref{def.ind.inclist}.

\begin{exercise}
\label{exe.Ialbe}Let $n\in\mathbb{N}$. Let $I$ be a subset of $\left\{
1,2,\ldots,n\right\}  $. Let $k=\left\vert I\right\vert $. Let $\left(
a_{1},a_{2},\ldots,a_{k}\right)  $ be the list of all elements of $I$ in
increasing order (with no repetitions). Let $\left(  b_{1},b_{2}%
,\ldots,b_{n-k}\right)  $ be the list of all elements of $\left\{
1,2,\ldots,n\right\}  \setminus I$ in increasing order (with no repetitions).
Let $\alpha\in S_{k}$ and $\beta\in S_{n-k}$. Prove the following:

\textbf{(a)} There exists a unique $\sigma\in S_{n}$ satisfying%
\[
\left(  \sigma\left(  1\right)  ,\sigma\left(  2\right)  ,\ldots,\sigma\left(
n\right)  \right)  =\left(  a_{\alpha\left(  1\right)  },a_{\alpha\left(
2\right)  },\ldots,a_{\alpha\left(  k\right)  },b_{\beta\left(  1\right)
},b_{\beta\left(  2\right)  },\ldots,b_{\beta\left(  n-k\right)  }\right)  .
\]
Denote this $\sigma$ by $\sigma_{I,\alpha,\beta}$.

\textbf{(b)} Let $\sum I$ denote the sum of all elements of $I$. (Thus, $\sum
I=\sum_{i\in I}i$.) We have%
\[
\ell\left(  \sigma_{I,\alpha,\beta}\right)  =\ell\left(  \alpha\right)
+\ell\left(  \beta\right)  +\sum I-\left(  1+2+\cdots+k\right)
\]
and%
\[
\left(  -1\right)  ^{\sigma_{I,\alpha,\beta}}=\left(  -1\right)  ^{\alpha
}\cdot\left(  -1\right)  ^{\beta}\cdot\left(  -1\right)  ^{\sum I-\left(
1+2+\cdots+k\right)  }.
\]

\textbf{(c)} Forget that we fixed $\alpha$ and $\beta$. We thus have defined
an element $\sigma_{I,\alpha,\beta}\in S_{n}$ for every $\alpha\in S_{k}$ and
every $\beta\in S_{n-k}$. The map%
\begin{align*}
S_{k}\times S_{n-k}  &  \rightarrow\left\{  \tau\in S_{n}\ \mid\ \tau\left(
\left\{  1,2,\ldots,k\right\}  \right)  =I\right\}  ,\\
\left(  \alpha,\beta\right)   &  \mapsto\sigma_{I,\alpha,\beta}%
\end{align*}
is well-defined and a bijection.
\end{exercise}

We can define transpositions not only in the symmetric group $S_{n}$, but also
more generally for arbitrary sets $X$:

\begin{definition}
\label{def.transposX}Let $X$ be a set. Let $i$ and $j$ be two distinct
elements of $X$. We let $t_{i,j}$ be the permutation of $X$ which swaps $i$
with $j$ while leaving all other elements of $X$ unchanged. Such a permutation
is called a \textit{transposition} of $X$.
\end{definition}

Clearly, Definition \ref{def.transposX} is a generalization of Definition
\ref{def.transpos}.

\begin{exercise}
\label{exe.perm.transX}Let $X$ be a finite set. Recall that if $\sigma$ is any
permutation of $X$, then the sign $\left(  -1\right)  ^{\sigma}$ of $\sigma$
is well-defined (by Exercise \ref{exe.ps4.2}). Prove the following:

\textbf{(a)} For any two distinct elements $i$ and $j$ of $X$, we have
$\left(  -1\right)  ^{t_{i,j}}=-1$.

\textbf{(b)} Any permutation of $X$ can be written as a composition of
finitely many transpositions of $X$.

\textbf{(c)} Let $\sigma$ be a permutation of $X$.\ Assume that $\sigma$ can
be written as a composition of $k$ transpositions of $X$. Then, $\left(
-1\right)  ^{\sigma}=\left(  -1\right)  ^{k}$.
\end{exercise}

\subsection{Cycles}

Next, we shall discuss another specific class of permutations: the
\textit{cycles}.

\begin{definition}
\label{def.perm.cycles}Let $n\in\mathbb{N}$. Let $\left[  n\right]  =\left\{
1,2,\ldots,n\right\}  $.

Let $k\in\left\{  1,2,\ldots,n\right\}  $. Let $i_{1},i_{2},\ldots,i_{k}$ be
$k$ distinct elements of $\left[  n\right]  $. We define $\operatorname*{cyc}%
\nolimits_{i_{1},i_{2},\ldots,i_{k}}$ to be the permutation in $S_{n}$ which
sends $i_{1},i_{2},\ldots,i_{k}$ to $i_{2},i_{3},\ldots,i_{k},i_{1}$,
respectively, while leaving all other elements of $\left[  n\right]  $ fixed.
In other words, we define $\operatorname*{cyc}\nolimits_{i_{1},i_{2}%
,\ldots,i_{k}}$ to be the permutation in $S_{n}$ given by%
\[
\left(
\begin{array}
[c]{r}%
\operatorname*{cyc}\nolimits_{i_{1},i_{2},\ldots,i_{k}}\left(  p\right)  =%
\begin{cases}
i_{j+1}, & \text{if }p=i_{j}\text{ for some }j\in\left\{  1,2,\ldots
,k\right\}  ;\\
p, & \text{otherwise}%
\end{cases}
\\
\ \ \ \ \ \ \ \ \ \ \text{for every }p\in\left[  n\right]
\end{array}
\right)  ,
\]
where $i_{k+1}$ means $i_{1}$.

(Again, the notation $\operatorname*{cyc}\nolimits_{i_{1},i_{2},\ldots,i_{k}}$
conceals the parameter $n$, which will hopefully not cause any confusion.)

A permutation of the form $\operatorname*{cyc}\nolimits_{i_{1},i_{2}%
,\ldots,i_{k}}$ is said to be a $k$\textit{-cycle} (or sometimes just a
\textit{cycle}, or a \textit{cyclic permutation}). Of course, the name stems
from the fact that it \textquotedblleft cycles\textquotedblright\ through the
elements $i_{1},i_{2},\ldots,i_{k}$ (by sending each of them to the next one
and the last one back to the first) and leaves all other elements unchanged.
\end{definition}

\begin{example}
\label{exa.perm.cycles}Let $n\in\mathbb{N}$. The following facts follow easily
from Definition \ref{def.perm.cycles}:

\textbf{(a)} For every $i\in\left\{  1,2,\ldots,n\right\}  $, we have
$\operatorname*{cyc}\nolimits_{i}=\operatorname*{id}$. In other words, any
$1$-cycle is the identity permutation $\operatorname*{id}$.

\textbf{(b)} If $i$ and $j$ are two distinct elements of $\left\{
1,2,\ldots,n\right\}  $, then $\operatorname*{cyc}\nolimits_{i,j}=t_{i,j}$.
(See Definition \ref{def.transpos} for the definition of $t_{i,j}$.)

\textbf{(c)} If $k\in\left\{  1,2,\ldots,n-1\right\}  $, then
$\operatorname*{cyc}\nolimits_{k,k+1}=s_{k}$.

\textbf{(d)} If $n=5$, then $\operatorname*{cyc}\nolimits_{2,5,3}$ is the
permutation which sends $1$ to $1$, $2$ to $5$, $3$ to $2$, $4$ to $4$, and
$5$ to $3$. (In other words, it is the permutation which is $\left(
1,5,2,4,3\right)  $ in one-line notation.)

\textbf{(e)} If $k\in\left\{  1,2,\ldots,n\right\}  $, and if $i_{1}%
,i_{2},\ldots,i_{k}$ are $k$ pairwise distinct elements of $\left[  n\right]
$, then%
\[
\operatorname*{cyc}\nolimits_{i_{1},i_{2},\ldots,i_{k}}=\operatorname*{cyc}%
\nolimits_{i_{2},i_{3},\ldots,i_{k},i_{1}}=\operatorname*{cyc}\nolimits_{i_{3}%
,i_{4},\ldots,i_{k},i_{1},i_{2}}=\cdots=\operatorname*{cyc}\nolimits_{i_{k}%
,i_{1},i_{2},\ldots,i_{k-1}}.
\]
(In less formal words: The $k$-cycle $\operatorname*{cyc}\nolimits_{i_{1}%
,i_{2},\ldots,i_{k}}$ does not change when we cyclically rotate the list
$\left(  i_{1},i_{2},\ldots,i_{k}\right)  $.)
\end{example}

\begin{remark}
What we called $\operatorname*{cyc}\nolimits_{i_{1},i_{2},\ldots,i_{k}}$ in
Definition \ref{def.perm.cycles} is commonly denoted by $\left(  i_{1}%
,i_{2},\ldots,i_{k}\right)  $ in the literature. But this latter notation
$\left(  i_{1},i_{2},\ldots,i_{k}\right)  $ would clash with one-line notation
for permutations (the cycle $\operatorname*{cyc}\nolimits_{1,2,3}\in S_{3}$ is
not the same as the permutation which is $\left(  1,2,3\right)  $ in one-line
notation) and also with the standard notation for $k$-tuples. This is why we
prefer to use the notation $\operatorname*{cyc}\nolimits_{i_{1},i_{2}%
,\ldots,i_{k}}$. (That said, we are not going to use $k$-cycles very often.)
\end{remark}

Any $k$-cycle is a composition of $k-1$ transpositions, as the following
exercise shows:

\begin{exercise}
\label{exe.perm.c=ttt}Let $n\in\mathbb{N}$. Let $\left[  n\right]  =\left\{
1,2,\ldots,n\right\}  $. Let $k\in\left\{  1,2,\ldots,n\right\}  $. Let
$i_{1},i_{2},\ldots,i_{k}$ be $k$ distinct elements of $\left[  n\right]  $.
Prove that%
\[
\operatorname*{cyc}\nolimits_{i_{1},i_{2},\ldots,i_{k}}=t_{i_{1},i_{2}}\circ
t_{i_{2},i_{3}}\circ\cdots\circ t_{i_{k-1},i_{k}}.
\]
(We are using Definition \ref{def.transpos} here.)
\end{exercise}

The following exercise gathers some further properties of cycles. Parts
\textbf{(a)} and \textbf{(d)} and, to a lesser extent, \textbf{(b)} are fairly
important and you should make sure you know how to solve them. The
significantly more difficult part \textbf{(c)} is more of a curiosity with an
interesting proof (I have not found an application of it so far; skip it if
you do not want to spend time on what is essentially a contest problem).

\begin{exercise}
\label{exe.perm.cycles}Let $n\in\mathbb{N}$. Let $\left[  n\right]  =\left\{
1,2,\ldots,n\right\}  $. Let $k\in\left\{  1,2,\ldots,n\right\}  $.

\textbf{(a)} For every $\sigma\in S_{n}$ and every $k$ distinct elements
$i_{1},i_{2},\ldots,i_{k}$ of $\left[  n\right]  $, prove that
\[
\sigma\circ\operatorname*{cyc}\nolimits_{i_{1},i_{2},\ldots,i_{k}}\circ
\sigma^{-1}=\operatorname*{cyc}\nolimits_{\sigma\left(  i_{1}\right)
,\sigma\left(  i_{2}\right)  ,\ldots,\sigma\left(  i_{k}\right)  }.
\]

\textbf{(b)} For every $p\in\left\{  0,1,\ldots,n-k\right\}  $, prove that%
\[
\ell\left(  \operatorname*{cyc}\nolimits_{p+1,p+2,\ldots,p+k}\right)  =k-1.
\]

\textbf{(c)} For every $k$ distinct elements $i_{1},i_{2},\ldots,i_{k}$ of
$\left[  n\right]  $, prove that
\[
\ell\left(  \operatorname*{cyc}\nolimits_{i_{1},i_{2},\ldots,i_{k}}\right)
\geq k-1.
\]

\textbf{(d)} For every $k$ distinct elements $i_{1},i_{2},\ldots,i_{k}$ of
$\left[  n\right]  $, prove that
\[
\left(  -1\right)  ^{\operatorname*{cyc}\nolimits_{i_{1},i_{2},\ldots,i_{k}}%
}=\left(  -1\right)  ^{k-1}.
\]

\end{exercise}

\begin{remark}
Exercise \ref{exe.perm.cycles} \textbf{(d)} shows that every $k$-cycle in
$S_{n}$ has sign $\left(  -1\right)  ^{k-1}$. However, the length of a
$k$-cycle need not be $k-1$. Exercise \ref{exe.perm.cycles} \textbf{(c)} shows
that this length is always $\geq k-1$, but it can take other values as well.
For instance, in $S_{4}$, the length of the $3$-cycle $\operatorname*{cyc}%
\nolimits_{1,4,3}$ is $4$. (Another example are the transpositions $t_{i,j}$
from Definition \ref{def.transpos}; these are $2$-cycles but can have length
$>1$.)

I don't know a simple way to describe when equality holds in Exercise
\ref{exe.perm.cycles} \textbf{(c)}. It holds whenever $i_{1},i_{2}%
,\ldots,i_{k}$ are consecutive integers (due to Exercise \ref{exe.perm.cycles}
\textbf{(b)}), but also in some other cases; for example, the $4$-cycle
$\operatorname*{cyc}\nolimits_{1,3,4,2}$ in $S_{4}$ has length $3$.
\end{remark}

\begin{remark}
\label{rmk.perm.cycles.decompose}The main reason why cycles are useful is
that, essentially, every permutation can be \textquotedblleft
decomposed\textquotedblright\ into cycles. We shall not use this fact, but
since it is generally important, let us briefly explain what it means. (You
will probably learn more about it in any standard course on abstract algebra.)

Fix $n\in\mathbb{N}$. Let $\left[  n\right]  =\left\{  1,2,\ldots,n\right\}
$. Two cycles $\alpha$ and $\beta$ in $S_{n}$ are said to be \textit{disjoint}
if they can be written as $\alpha=\operatorname*{cyc}\nolimits_{i_{1}%
,i_{2},\ldots,i_{k}}$ and $\beta=\operatorname*{cyc}\nolimits_{j_{1}%
,j_{2},\ldots,j_{\ell}}$ for $k+\ell$ distinct elements $i_{1},i_{2}%
,\ldots,i_{k},j_{1},j_{2},\ldots,j_{\ell}$ of $\left[  n\right]  $. For
example, the two cycles $\operatorname*{cyc}\nolimits_{1,3}$ and
$\operatorname*{cyc}\nolimits_{2,6,7}$ in $S_{8}$ are disjoint, but the two
cycles $\operatorname*{cyc}\nolimits_{1,4}$ and $\operatorname*{cyc}%
\nolimits_{2,4}$ are not. It is easy to see that any two disjoint cycles
$\alpha$ and $\beta$ commute (i.e., satisfy $\alpha\circ\beta=\beta\circ
\alpha$). Therefore, when you see a composition $\alpha_{1}\circ\alpha
_{2}\circ\cdots\circ\alpha_{p}$ of several pairwise disjoint cycles, you can
reorder its factors arbitrarily without changing the result (for example,
$\alpha_{3}\circ\alpha_{1}\circ\alpha_{4}\circ\alpha_{2}=\alpha_{1}\circ
\alpha_{2}\circ\alpha_{3}\circ\alpha_{4}$ if $p=4$).

Now, the fact I am talking about says the following: Every permutation in
$S_{n}$ can be written as a composition of several pairwise disjoint cycles.
For example, let $n=9$, and let $\sigma\in S_{9}$ be the permutation which is
written $\left(  4,6,1,3,5,2,9,8,7\right)  $ in one-line notation (i.e., we
have $\sigma\left(  1\right)  =4$, $\sigma\left(  2\right)  =6$, etc.). Then,
$\sigma$ can be written as a composition of several pairwise disjoint cycles
as follows:%
\begin{equation}
\sigma=\operatorname*{cyc}\nolimits_{1,4,3}\circ\operatorname*{cyc}%
\nolimits_{7,9}\circ\operatorname*{cyc}\nolimits_{2,6}.
\label{eq.rmk.perm.cycles.decompose.1}%
\end{equation}
Indeed, here is how such a decomposition can be found: Let us draw a directed
graph whose vertices are $1,2,\ldots,n$, and which has an arc $i\rightarrow
\sigma\left(  i\right)  $ for every $i\in\left[  n\right]  $. (Thus, it has
$n$ arcs altogether; some of them can be loops.) For our permutation
$\sigma\in S_{9}$, this graph looks as follows:
\begin{equation}
\begin{tikzpicture}[->,>=stealth',shorten >=1pt,auto,node distance=3cm, thick,main node/.style={circle,fill=blue!20,draw}] \node[main node] (1) {1}; \node[main node] (3) [below right of=1] {3}; \node[main node] (4) [above right of=1] {4}; \node[main node] (2) [above left of=1] {2}; \node[main node] (6) [below left of=1] {6}; \node[main node] (5) [below right of=4] {5}; \node[main node] (7) [above right of=5] {7}; \node[main node] (9) [below right of=5] {9}; \node[main node] (8) [above right of=9] {8}; \path[every node/.style={font=\sffamily\small}] (1) edge (4) (2) edge [bend right] (6) (3) edge (1) (4) edge (3) (5) edge [loop left] (5) (6) edge [bend right] (2) (7) edge [bend right] (9) (8) edge [loop left] (8) (9) edge [bend right] (7); \end{tikzpicture}.
\label{eq.rmk.perm.cycles.decompose.2}%
\end{equation}
Obviously, at each vertex $i$ of this graph, exactly one arc begins (namely,
the arc $i\rightarrow\sigma\left(  i\right)  $). Moreover, since $\sigma$ is
invertible, it is also clear that at each vertex $i$ of this graph, exactly
one arc ends (namely, the arc $\sigma^{-1}\left(  i\right)  \rightarrow i$).
Due to the way we constructed this graph, it is clear that it completely
describes our permutation $\sigma$: Namely, if we want to find $\sigma\left(
i\right)  $ for a given $i\in\left[  n\right]  $, we should just locate the
vertex $i$ on the graph, and follow the arc that begins at this vertex; the
endpoint of this arc will be $\sigma\left(  i\right)  $.

Now, a look at this graph reveals five directed cycles (in the sense of
\textquotedblleft paths which end at the same vertex at which they
begin\textquotedblright, not yet in the sense of \textquotedblleft cyclic
permutations\textquotedblright). The first one passes through the vertices $2$
and $6$; the second passes through the vertices $3$, $1$ and $4$; the third,
through the vertex $5$ (it is what is called a \textquotedblleft trivial
cycle\textquotedblright), and so on. To each of these cycles we can assign a
cyclic permutation in $S_{n}$: namely, if the cycle passes through the
vertices $i_{1},i_{2},\ldots,i_{k}$ (in this order, and with no repetitions),
then we assign to it the cyclic permutation $\operatorname*{cyc}%
\nolimits_{i_{1},i_{2},\ldots,i_{k}}\in S_{n}$. The cyclic permutations
assigned to all five directed cycles are pairwise disjoint, and their
composition is%
\[
\operatorname*{cyc}\nolimits_{2,6}\circ\operatorname*{cyc}\nolimits_{3,1,4}%
\circ\operatorname*{cyc}\nolimits_{5}\circ\operatorname*{cyc}\nolimits_{7,9}%
\circ\operatorname*{cyc}\nolimits_{8}.
\]
But this composition must be $\sigma$ (because if we apply this composition to
an element $i\in\left[  n\right]  $, then we obtain the \textquotedblleft next
vertex after $i$\textquotedblright\ on the directed cycle which passes through
$i$; but due to how we constructed our graph, this \textquotedblleft next
vertex\textquotedblright\ will be precisely $\sigma\left(  i\right)  $).
Hence, we have%
\begin{equation}
\sigma=\operatorname*{cyc}\nolimits_{2,6}\circ\operatorname*{cyc}%
\nolimits_{3,1,4}\circ\operatorname*{cyc}\nolimits_{5}\circ\operatorname*{cyc}%
\nolimits_{7,9}\circ\operatorname*{cyc}\nolimits_{8}.
\label{eq.rmk.perm.cycles.decompose.3}%
\end{equation}
Thus, we have found a way to write $\sigma$ as a composition of several
pairwise disjoint cycles. We can rewrite (and even simplify) this
representation a bit: Namely, we can simplify
(\ref{eq.rmk.perm.cycles.decompose.3}) by removing the factors
$\operatorname*{cyc}\nolimits_{5}$ and $\operatorname*{cyc}\nolimits_{8}$
(because both of these factors equal $\operatorname*{id}$); thus we obtain
$\sigma=\operatorname*{cyc}\nolimits_{2,6}\circ\operatorname*{cyc}%
\nolimits_{3,1,4}\circ\operatorname*{cyc}\nolimits_{7,9}$. We can furthermore
swap $\operatorname*{cyc}\nolimits_{2,6}$ with $\operatorname*{cyc}%
\nolimits_{3,1,4}$ (since disjoint cycles commute), therefore obtaining
$\sigma=\operatorname*{cyc}\nolimits_{3,1,4}\circ\operatorname*{cyc}%
\nolimits_{2,6}\circ\operatorname*{cyc}\nolimits_{7,9}$. Next, we can swap
$\operatorname*{cyc}\nolimits_{2,6}$ with $\operatorname*{cyc}\nolimits_{7,9}%
$, obtaining $\sigma=\operatorname*{cyc}\nolimits_{3,1,4}\circ
\operatorname*{cyc}\nolimits_{7,9}\circ\operatorname*{cyc}\nolimits_{2,6}$.
Finally, we can rewrite $\operatorname*{cyc}\nolimits_{3,1,4}$ as
$\operatorname*{cyc}\nolimits_{1,4,3}$, and we obtain
(\ref{eq.rmk.perm.cycles.decompose.1}).

In general, for every $n\in\mathbb{N}$, every permutation $\sigma\in S_{n}$
can be represented as a composition of several pairwise disjoint cycles (which
can be found by drawing a directed graph as in our example above). This
representation is not literally unique, because we can modify it by:

\begin{itemize}
\item adding or removing trivial factors (i.e., factors of the form
$\operatorname*{cyc}\nolimits_{i}=\operatorname*{id}$);

\item swapping different cycles;

\item rewriting $\operatorname*{cyc}\nolimits_{i_{1},i_{2},\ldots,i_{k}}$ as
$\operatorname*{cyc}\nolimits_{i_{2},i_{3},\ldots,i_{k},i_{1}}$.
\end{itemize}

However, it is unique \textbf{up to these modifications}; in other words, any
two representations of $\sigma$ as a composition of several pairwise disjoint
cycles can be transformed into one another by such modifications.

The proofs of all these statements are fairly easy. (One does have to check
certain things, e.g., that the directed graph really consists of disjoint
directed cycles. For a complete proof, see \cite[Theorem 1.5.3]{Goodman} or
\cite[Chapter I, \S 5.7, Proposition 7]{Bourba74} or \cite[\S 1.9, proof of
Theorem 1.5.1]{Sagan19} or various other texts on algebra.)

Representing a permutation $\sigma\in S_{n}$ as a composition of several
pairwise disjoint cycles can be done very quickly, and thus gives a quick way
to find $\left(  -1\right)  ^{\sigma}$ (because Exercise \ref{exe.perm.cycles}
\textbf{(d)} tells us how to find the sign of a $k$-cycle). This is
significantly faster than counting inversions of $\sigma$.
\end{remark}

\subsection{\label{sect.perm.lehmer}The Lehmer code}

In this short section, we shall introduce the \textit{Lehmer code} of a
permutation. Throughout Section \ref{sect.perm.lehmer}, we will use the
following notations:

\begin{definition}
\textbf{(a)} Whenever $m$ is an integer, we shall use the notation $\left[
m\right]  $ for the set $\left\{  1,2,\ldots,m\right\}  $. (This is an empty
set when $m\leq0$.)

\textbf{(b)} Whenever $m$ is an integer, we shall use the notation $\left[
m\right]  _{0}$ for the set $\left\{  0,1,\ldots,m\right\}  $. (This is an
empty set when $m<0$.)
\end{definition}

\begin{definition}
Let $n\in\mathbb{N}$. We consider $n$ to be fixed throughout Section
\ref{sect.perm.lehmer}.

Let $H$ denote the set $\left[  n-1\right]  _{0}\times\left[  n-2\right]
_{0}\times\cdots\times\left[  n-n\right]  _{0}$.
\end{definition}

\begin{definition}
Let $\sigma\in S_{n}$ and $i\in\left[  n\right]  $. Then, $\ell_{i}\left(
\sigma\right)  $ shall denote the number of all $j\in\left\{  i+1,i+2,\ldots
,n\right\}  $ such that $\sigma\left(  i\right)  >\sigma\left(  j\right)  $.
\end{definition}

\begin{example}
For this example, set $n=5$, and let $\sigma\in S_{5}$ be the permutation that
sends $1,2,3,4,5$ to $4,3,2,1,5$. Then, $\ell_{2}\left(  \sigma\right)  $ is
the number of all $j\in\left\{  3,4,5\right\}  $ such that $\sigma\left(
2\right)  >\sigma\left(  j\right)  $. These $j$ are $3$ and $4$ (because
$\sigma\left(  2\right)  >\sigma\left(  3\right)  $ and $\sigma\left(
2\right)  >\sigma\left(  4\right)  $ but not $\sigma\left(  2\right)
>\sigma\left(  5\right)  $); therefore, $\ell_{2}\left(  \sigma\right)  =2$.
Similarly, $\ell_{1}\left(  \sigma\right)  =3$, $\ell_{3}\left(
\sigma\right)  =1$, $\ell_{4}\left(  \sigma\right)  =0$ and $\ell_{5}\left(
\sigma\right)  =0$.
\end{example}

The following two facts are almost trivial:\footnote{See Exercise
\ref{exe.perm.lehmer.prove} below for the proofs of all the following
results.}

\begin{proposition}
\label{prop.perm.lehmer.l}Let $\sigma\in S_{n}$. Then, $\ell\left(
\sigma\right)  =\ell_{1}\left(  \sigma\right)  +\ell_{2}\left(  \sigma\right)
+\cdots+\ell_{n}\left(  \sigma\right)  $.
\end{proposition}

\begin{proposition}
\label{prop.perm.lehmer.wd}Let $\sigma\in S_{n}$. Then, $\left(  \ell
_{1}\left(  \sigma\right)  ,\ell_{2}\left(  \sigma\right)  ,\ldots,\ell
_{n}\left(  \sigma\right)  \right)  \in H$.
\end{proposition}

The following simple lemma gives two equivalent ways to define $\ell
_{i}\left(  \sigma\right)  $ for $\sigma\in S_{n}$ and $i\in\left[  n\right]
$:

\begin{lemma}
\label{lem.perm.lexico1.lis}Let $\sigma\in S_{n}$ and $i\in\left[  n\right]
$. Then:

\textbf{(a)} We have $\ell_{i}\left(  \sigma\right)  =\left\vert \left[
\sigma\left(  i\right)  -1\right]  \setminus\sigma\left(  \left[  i\right]
\right)  \right\vert $.

\textbf{(b)} We have $\ell_{i}\left(  \sigma\right)  =\left\vert \left[
\sigma\left(  i\right)  -1\right]  \setminus\sigma\left(  \left[  i-1\right]
\right)  \right\vert $.

\textbf{(c)} We have $\sigma\left(  i\right)  \leq i+\ell_{i}\left(
\sigma\right)  $.
\end{lemma}

Before we state the next proposition, we introduce another notation:

\begin{definition}
\label{def.perm.lehmer.lex-ord}Let $\left(  a_{1},a_{2},\ldots,a_{n}\right)  $
and $\left(  b_{1},b_{2},\ldots,b_{n}\right)  $ be two $n$-tuples of integers.
We say that $\left(  a_{1},a_{2},\ldots,a_{n}\right)  <_{\operatorname*{lex}%
}\left(  b_{1},b_{2},\ldots,b_{n}\right)  $ if and only if there exists some
$k\in\left[  n\right]  $ such that $a_{k}\neq b_{k}$, and the
\textbf{smallest} such $k$ satisfies $a_{k}<b_{k}$.
\end{definition}

For example, $\left(  4,1,2,5\right)  <_{\operatorname*{lex}}\left(
4,1,3,0\right)  $ and $\left(  1,1,0,1\right)  <_{\operatorname*{lex}}\left(
2,0,0,0\right)  $. The relation $<_{\operatorname*{lex}}$ is usually
pronounced \textquotedblleft is lexicographically smaller
than\textquotedblright; the word \textquotedblleft
lexicographic\textquotedblright\ comes from the idea that if numbers were
letters, then a \textquotedblleft word\textquotedblright\ $a_{1}a_{2}\cdots
a_{n}$ would appear earlier in a dictionary than $b_{1}b_{2}\cdots b_{n}$ if
and only if $\left(  a_{1},a_{2},\ldots,a_{n}\right)  <_{\operatorname*{lex}%
}\left(  b_{1},b_{2},\ldots,b_{n}\right)  $.

\begin{proposition}
\label{prop.perm.lehmer.lex}Let $\sigma\in S_{n}$ and $\tau\in S_{n}$ be such
that
\[
\left(  \sigma\left(  1\right)  ,\sigma\left(  2\right)  ,\ldots,\sigma\left(
n\right)  \right)  <_{\operatorname*{lex}}\left(  \tau\left(  1\right)
,\tau\left(  2\right)  ,\ldots,\tau\left(  n\right)  \right)  .
\]
Then,
\[
\left(  \ell_{1}\left(  \sigma\right)  ,\ell_{2}\left(  \sigma\right)
,\ldots,\ell_{n}\left(  \sigma\right)  \right)  <_{\operatorname*{lex}}\left(
\ell_{1}\left(  \tau\right)  ,\ell_{2}\left(  \tau\right)  ,\ldots,\ell
_{n}\left(  \tau\right)  \right)  .
\]

\end{proposition}

We can now define the Lehmer code:

\begin{definition}
\label{def.perm.lehmer.L}Define the map $L:S_{n}\rightarrow H$ by%
\[
\left(  L\left(  \sigma\right)  =\left(  \ell_{1}\left(  \sigma\right)
,\ell_{2}\left(  \sigma\right)  ,\ldots,\ell_{n}\left(  \sigma\right)
\right)  \ \ \ \ \ \ \ \ \ \ \text{for each }\sigma\in S_{n}\right)  .
\]
(This is well-defined because of Proposition \ref{prop.perm.lehmer.wd}.)

If $\sigma\in S_{n}$ is any permutation, then $L\left(  \sigma\right)
=\left(  \ell_{1}\left(  \sigma\right)  ,\ell_{2}\left(  \sigma\right)
,\ldots,\ell_{n}\left(  \sigma\right)  \right)  $ is called the \textit{Lehmer
code} of $\sigma$.
\end{definition}

\begin{theorem}
\label{thm.perm.lehmer.bij}The map $L:S_{n}\rightarrow H$ is a bijection.
\end{theorem}

Using this theorem and Proposition \ref{prop.perm.lehmer.l}, we can easily
show the following:

\begin{corollary}
\label{cor.perm.lehmer.lensum}We have%
\[
\sum_{w\in S_{n}}x^{\ell\left(  w\right)  }=\left(  1+x\right)  \left(
1+x+x^{2}\right)  \cdots\left(  1+x+x^{2}+\cdots+x^{n-1}\right)
\]
(an equality between polynomials in $x$). (The right hand side of this
equality should be understood as the empty product when $n\leq1$.)
\end{corollary}

\begin{exercise}
\label{exe.perm.lehmer.prove}Prove Proposition \ref{prop.perm.lehmer.l},
Proposition \ref{prop.perm.lehmer.wd}, Lemma \ref{lem.perm.lexico1.lis},
Proposition \ref{prop.perm.lehmer.lex}, Theorem \ref{thm.perm.lehmer.bij} and
Corollary \ref{cor.perm.lehmer.lensum}.
\end{exercise}

See \cite[\S 2.1]{Manive01} and \cite[\S 11.3]{Kerber99} for further
properties of permutations related to the Lehmer code. (In particular,
\cite[proof of Proposition 2.1.2]{Manive01} and \cite[Corollary 11.3.5]%
{Kerber99} give two different ways of reconstructing a permutation from its
Lehmer code; moreover, \cite[Corollary 11.3.5]{Kerber99} shows how the Lehmer
code of a permutation $\sigma\in S_{n}$ leads to a specific representation of
$\sigma$ as a product of some of the $s_{1},s_{2},\ldots,s_{n-1}$.)

\begin{exercise}
\label{exe.perm.lisitau}Let $n\in\mathbb{N}$. Let $\sigma\in S_{n}$ and
$\tau\in S_{n}$. We shall use the notation from Definition \ref{def.iverson}.

\textbf{(a)} Prove that each $i\in\left[  n\right]  $ satisfies%
\begin{align*}
&  \ell_{\tau\left(  i\right)  }\left(  \sigma\right)  +\ell_{i}\left(
\tau\right)  -\ell_{i}\left(  \sigma\circ\tau\right) \\
&  =\sum_{j\in\left[  n\right]  }\left[  j>i\right]  \left[  \tau\left(
i\right)  >\tau\left(  j\right)  \right]  \left[  \sigma\left(  \tau\left(
j\right)  \right)  >\sigma\left(  \tau\left(  i\right)  \right)  \right] \\
&  \ \ \ \ \ \ \ \ \ \ +\sum_{j\in\left[  n\right]  }\left[  i>j\right]
\left[  \tau\left(  j\right)  >\tau\left(  i\right)  \right]  \left[
\sigma\left(  \tau\left(  i\right)  \right)  >\sigma\left(  \tau\left(
j\right)  \right)  \right]  .
\end{align*}

\textbf{(b)} Prove that%
\[
\ell\left(  \sigma\right)  +\ell\left(  \tau\right)  -\ell\left(  \sigma
\circ\tau\right)  =2\sum_{i\in\left[  n\right]  }\sum_{j\in\left[  n\right]
}\left[  j>i\right]  \left[  \tau\left(  i\right)  >\tau\left(  j\right)
\right]  \left[  \sigma\left(  \tau\left(  j\right)  \right)  >\sigma\left(
\tau\left(  i\right)  \right)  \right]  .
\]

\textbf{(c)} Give a new solution to Exercise \ref{exe.ps2.2.5} \textbf{(a)}.

\textbf{(d)} Give a new solution to Exercise \ref{exe.ps2.2.5} \textbf{(b)}.

\textbf{(e)} Give a new solution to Exercise \ref{exe.ps2.2.5} \textbf{(c)}.
\end{exercise}

\begin{exercise}
\label{exe.perm.lisitij}Let $n\in\mathbb{N}$ and $\sigma\in S_{n}$. Let $i$
and $j$ be two elements of $\left[  n\right]  $ such that $i<j$ and
$\sigma\left(  i\right)  >\sigma\left(  j\right)  $. Let $Q$ be the set of all
$k\in\left\{  i+1,i+2,\ldots,j-1\right\}  $ satisfying $\sigma\left(
i\right)  >\sigma\left(  k\right)  >\sigma\left(  j\right)  $. Prove that
\[
\ell\left(  \sigma\circ t_{i,j}\right)  =\ell\left(  \sigma\right)
-2\left\vert Q\right\vert -1.
\]

\end{exercise}

The following exercise shows an explicit way of expressing every permutation
$\sigma\in S_{n}$ as a product of $\ell\left(  \sigma\right)  $ many simple
transpositions (i.e., transpositions of the form $s_{i}$ with $i\in\left\{
1,2,\ldots,n-1\right\}  $):

\begin{exercise}
\label{exe.perm.lehmer.rothe}Let $n\in\mathbb{N}$ and $\sigma\in S_{n}$. If
$u$ and $v$ are any two elements of $\left[  n\right]  $ such that $u\leq v$,
then we define a permutation $c_{u,v}\in S_{n}$ by%
\[
c_{u,v}=\operatorname*{cyc}\nolimits_{v,v-1,v-2,\ldots,u}.
\]
For each $i\in\left[  n\right]  $, we define a permutation $a_{i}\in S_{n}$
by
\[
a_{i}=c_{i,i+\ell_{i}\left(  \sigma\right)  }.
\]

\textbf{(a)} Prove that $a_{i}$ is well-defined for each $i\in\left[
n\right]  $.

\textbf{(b)} Prove that each $i\in\left[  n\right]  $ satisfies $a_{i}%
=s_{i^{\prime}-1}\circ s_{i^{\prime}-2}\circ\cdots\circ s_{i}$, where
$i^{\prime}=i+\ell_{i}\left(  \sigma\right)  $.

\textbf{(c)} Prove that $\sigma=a_{1}\circ a_{2}\circ\cdots\circ a_{n}$.

\textbf{(d)} Solve Exercise \ref{exe.ps2.2.5} \textbf{(e)} again.

\textbf{(e)} Solve Exercise \ref{exe.ps2.2.4} \textbf{(c)} again.
\end{exercise}

\subsection{Extending permutations}

In this short section, we shall discuss a simple yet useful concept: that of
extending a permutation of a set $Y$ to a larger set $X$ (where
\textquotedblleft larger\textquotedblright\ means that $Y\subseteq X$). The
following notations will be used throughout this section:

\begin{definition}
\label{def.perm.extend.SX}Let $X$ be a set. Then, $S_{X}$ denotes the set of
all permutations of $X$.
\end{definition}

\begin{definition}
\label{def.perm.extend.YtX}Let $X$ be a set. Let $Y$ be a subset of $X$. For
every map $\sigma:Y\rightarrow Y$, we define a map $\sigma^{\left(
Y\rightarrow X\right)  }:X\rightarrow X$ by%
\[
\left(  \sigma^{\left(  Y\rightarrow X\right)  }\left(  x\right)  =%
\begin{cases}
\sigma\left(  x\right)  , & \text{if }x\in Y;\\
x, & \text{if }x\notin Y
\end{cases}
\ \ \ \ \ \ \ \ \ \ \text{for every }x\in X\right)  .
\]
(This map $\sigma^{\left(  Y\rightarrow X\right)  }$ is indeed well-defined,
according to Proposition \ref{prop.perm.extend.YtX.wd} below.)
\end{definition}

The latter of these two definitions relies on the following lemma:

\begin{proposition}
\label{prop.perm.extend.YtX.wd}Let $X$ be a set. Let $Y$ be a subset of $X$.
Let $\sigma:Y\rightarrow Y$ be a map. Then, the map $\sigma^{\left(
Y\rightarrow X\right)  }$ in Definition \ref{def.perm.extend.YtX} is well-defined.
\end{proposition}

Proposition \ref{prop.perm.extend.YtX.wd} is easy to prove; its proof is part
of Exercise \ref{exe.perm.extend.proofs} further below.

The idea behind the definition of $\sigma^{\left(  Y\rightarrow X\right)  }$
in Definition \ref{def.perm.extend.YtX} is simple: $\sigma^{\left(
Y\rightarrow X\right)  }$ is just the most straightforward way of extending
$\sigma:Y\rightarrow Y$ to a map from $X$ to $X$ (namely, by letting it keep
every element of $X\setminus Y$ unchanged).

\begin{example}
\textbf{(a)} If $X=\left\{  1,2,3,4,5,6,7\right\}  $ and $Y=\left\{
1,2,3,4\right\}  $, and if $\sigma\in S_{Y}=S_{4}$ is the permutation whose
one-line notation is $\left(  4,1,3,2\right)  $, then $\sigma^{\left(
Y\rightarrow X\right)  }\in S_{X}=S_{7}$ is the permutation whose one-line
notation is $\left(  4,1,3,2,5,6,7\right)  $.

\textbf{(b)} More generally, if $X=\left\{  1,2,\ldots,n\right\}  $ and
$Y=\left\{  1,2,\ldots,m\right\}  $ for two nonnegative integers $n$ and $m$
satisfying $n\geq m$, and if $\sigma\in S_{Y}=S_{m}$ is any permutation, then
the permutation $\sigma^{\left(  Y\rightarrow X\right)  }\in S_{X}=S_{n}$ has
one-line notation $\left(  \sigma\left(  1\right)  ,\sigma\left(  2\right)
,\ldots,\sigma\left(  m\right)  ,m+1,m+2,\ldots,n\right)  $.

\textbf{(c)} If $X=\left\{  1,2,3,\ldots\right\}  $ and $Y=\left\{
1,2,\ldots,n\right\}  $ for some $n\in\mathbb{N}$, and if $\sigma\in
S_{Y}=S_{n}$ is any permutation, then the permutation $\sigma^{\left(
Y\rightarrow X\right)  }\in S_{X}=S_{\infty}$ is precisely the permutation
$\sigma_{\left(  \infty\right)  }$ defined in Remark
\ref{rmk.perm.inf.lazy-ext}.
\end{example}

Here are some further properties of the operation that transforms $\sigma$
into $\sigma^{\left(  Y\rightarrow X\right)  }$:

\begin{proposition}
\label{prop.perm.extend.YtX.inj-sur}Let $X$ be a set. Let $Y$ be a subset of
$X$.

\textbf{(a)} If $\alpha:Y\rightarrow Y$ and $\beta:Y\rightarrow Y$ are two
maps, then%
\[
\left(  \alpha\circ\beta\right)  ^{\left(  Y\rightarrow X\right)  }%
=\alpha^{\left(  Y\rightarrow X\right)  }\circ\beta^{\left(  Y\rightarrow
X\right)  }.
\]

\textbf{(b)} The map $\operatorname*{id}\nolimits_{Y}:Y\rightarrow Y$
satisfies%
\[
\left(  \operatorname*{id}\nolimits_{Y}\right)  ^{\left(  Y\rightarrow
X\right)  }=\operatorname*{id}\nolimits_{X}.
\]

\textbf{(c)} Every permutation $\sigma\in S_{Y}$ satisfies $\sigma^{\left(
Y\rightarrow X\right)  }\in S_{X}$ and%
\[
\left(  \sigma^{-1}\right)  ^{\left(  Y\rightarrow X\right)  }=\left(
\sigma^{\left(  Y\rightarrow X\right)  }\right)  ^{-1}.
\]

\textbf{(d)} We have%
\[
\left\{  \delta^{\left(  Y\rightarrow X\right)  }\ \mid\ \delta\in
S_{Y}\right\}  =\left\{  \tau\in S_{X}\ \mid\ \tau\left(  z\right)  =z\text{
for every }z\in X\setminus Y\right\}  .
\]

\textbf{(e)} The map%
\begin{align*}
S_{Y}  &  \rightarrow\left\{  \tau\in S_{X}\ \mid\ \tau\left(  z\right)
=z\text{ for every }z\in X\setminus Y\right\}  ,\\
\delta &  \mapsto\delta^{\left(  Y\rightarrow X\right)  }%
\end{align*}
is well-defined and bijective.
\end{proposition}

\begin{proposition}
\label{prop.perm.extend.YtX.ZtY}Let $X$ be a set. Let $Y$ be a subset of $X$.
Let $Z$ be a subset of $Y$. Let $\sigma:Z\rightarrow Z$ be any map. Then,%
\[
\left(  \sigma^{\left(  Z\rightarrow Y\right)  }\right)  ^{\left(
Y\rightarrow X\right)  }=\sigma^{\left(  Z\rightarrow X\right)  }.
\]

\end{proposition}

\begin{proposition}
\label{prop.perm.extend.YtX.commute}Let $X$ be a set. Let $Y$ be a subset of
$X$. Let $\alpha:Y\rightarrow Y$ be a map. Let $\beta:X\setminus Y\rightarrow
X\setminus Y$ be a map. Then,%
\[
\alpha^{\left(  Y\rightarrow X\right)  }\circ\beta^{\left(  X\setminus
Y\rightarrow X\right)  }=\beta^{\left(  X\setminus Y\rightarrow X\right)
}\circ\alpha^{\left(  Y\rightarrow X\right)  }.
\]

\end{proposition}

The above propositions are fairly straightforward; again, see Exercise
\ref{exe.perm.extend.proofs} for their proofs. Interestingly, we can use these
simple facts to prove the following nontrivial theorem:

\begin{theorem}
\label{thm.perm.extend.conj-inv.2}Let $X$ be a finite set. Let $\pi\in S_{X}$.
Then, there exists a $\sigma\in S_{X}$ such that $\sigma\circ\pi\circ
\sigma^{-1}=\pi^{-1}$.
\end{theorem}

In fact, we can show (by induction on $\left\vert X\right\vert $) the
following more general fact:

\begin{proposition}
\label{prop.perm.extend.conj-inv.1}Let $X$ be a finite set. Let $x\in X$. Let
$\pi\in S_{X}$. Then, there exists a $\sigma\in\left\{  \delta^{\left(
X\setminus\left\{  x\right\}  \rightarrow X\right)  }\ \mid\ \delta\in
S_{X\setminus\left\{  x\right\}  }\right\}  $ such that $\sigma\circ\pi
\circ\sigma^{-1}=\pi^{-1}$.
\end{proposition}

\begin{exercise}
\label{exe.perm.extend.proofs}Prove Proposition \ref{prop.perm.extend.YtX.wd},
Proposition \ref{prop.perm.extend.YtX.inj-sur}, Proposition
\ref{prop.perm.extend.YtX.ZtY}, Proposition \ref{prop.perm.extend.YtX.commute}%
, Proposition \ref{prop.perm.extend.conj-inv.1} and Theorem
\ref{thm.perm.extend.conj-inv.2}.
\end{exercise}

Theorem \ref{thm.perm.extend.conj-inv.2} is a known fact, and it is commonly
obtained as part of the study of conjugacy in symmetric groups. If $\pi_{1}$
and $\pi_{2}$ are two permutations of a set $X$, then $\pi_{1}$ is said to be
\textit{conjugate} to $\pi_{2}$ if and only if there exists some $\sigma\in
S_{X}$ such that $\sigma\circ\pi_{1}\circ\sigma^{-1}=\pi_{2}$. Thus, Theorem
\ref{thm.perm.extend.conj-inv.2} says that every permutation $\pi$ of a finite
set $X$ is conjugate to its inverse $\pi^{-1}$. Standard proofs of this
theorem\footnote{which, incidentally, also holds for infinite sets $X$,
provided that one believes in the Axiom of Choice} tend to derive it from the
fact that two permutations $\pi_{1}$ and $\pi_{2}$ of a finite set $X$ are
conjugate to one another\footnote{It is easy to see that being conjugate is a
symmetric relation: If a permutation $\pi_{1}$ is conjugate to a permutation
$\pi_{2}$, then $\pi_{2}$ is, in turn, conjugate to $\pi_{1}$.} if and only if
they have the same \textquotedblleft cycle type\textquotedblright\ (see
\cite[Theorem 5.7]{C-conjclass} for what this means and for a proof).

\subsection{Additional exercises}

Permutations and symmetric groups are a staple of combinatorics; there are
countless results involving them. For an example, B\'{o}na's book
\cite{Bona22}, as well as significant parts of Stanley's \cite{Stanley-EC1}
and \cite{Stanley-EC2} are devoted to them. In this section, I shall only give
a haphazard selection of exercises, which are not relevant to the rest of
these notes (thus can be skipped at will). I am not planning to provide
solutions for all of them.

\begin{exercise}
\label{exe.perm.aj-ai-sum}Let $n\in\mathbb{N}$. Let $\sigma\in S_{n}$. Let
$a_{1},a_{2},\ldots,a_{n}$ be any $n$ numbers. (Here, \textquotedblleft
number\textquotedblright\ means \textquotedblleft real
number\textquotedblright\ or \textquotedblleft complex
number\textquotedblright\ or \textquotedblleft rational
number\textquotedblright, as you prefer; this makes no difference.) Prove that%
\[
\sum_{\substack{1\leq i<j\leq n;\\\sigma\left(  i\right)  >\sigma\left(
j\right)  }}\left(  a_{j}-a_{i}\right)  =\sum_{i=1}^{n}a_{i}\left(
i-\sigma\left(  i\right)  \right)  .
\]
[Here, the summation sign \textquotedblleft$\sum_{\substack{1\leq i<j\leq
n;\\\sigma\left(  i\right)  >\sigma\left(  j\right)  }}$\textquotedblright%
\ means \textquotedblleft$\sum_{\substack{\left(  i,j\right)  \in\left\{
1,2,\ldots,n\right\}  ^{2};\\i<j\text{ and }\sigma\left(  i\right)
>\sigma\left(  j\right)  }}$\textquotedblright; this is a sum over all
inversions of $\sigma$.]
\end{exercise}

\begin{exercise}
\label{exe.perm.pij-pii-sum}Let $n\in\mathbb{N}$. Let $\pi\in S_{n}$.

\textbf{(a)} Prove that
\[
\sum_{\substack{1\leq i<j\leq n;\\\pi\left(  i\right)  >\pi\left(  j\right)
}}\left(  \pi\left(  j\right)  -\pi\left(  i\right)  \right)  =\sum
_{\substack{1\leq i<j\leq n;\\\pi\left(  i\right)  >\pi\left(  j\right)
}}\left(  i-j\right)  .
\]
[Here, the summation sign \textquotedblleft$\sum_{\substack{1\leq i<j\leq
n;\\\pi\left(  i\right)  >\pi\left(  j\right)  }}$\textquotedblright\ means
\textquotedblleft$\sum_{\substack{\left(  i,j\right)  \in\left\{
1,2,\ldots,n\right\}  ^{2};\\i<j\text{ and }\pi\left(  i\right)  >\pi\left(
j\right)  }}$\textquotedblright; this is a sum over all inversions of $\pi$.]

\textbf{(b)} Prove that%
\[
\sum_{\substack{1\leq i<j\leq n;\\\pi\left(  i\right)  <\pi\left(  j\right)
}}\left(  \pi\left(  j\right)  -\pi\left(  i\right)  \right)  =\sum
_{\substack{1\leq i<j\leq n;\\\pi\left(  i\right)  <\pi\left(  j\right)
}}\left(  j-i\right)  .
\]
[Here, the summation sign \textquotedblleft$\sum_{\substack{1\leq i<j\leq
n;\\\pi\left(  i\right)  <\pi\left(  j\right)  }}$\textquotedblright\ means
\textquotedblleft$\sum_{\substack{\left(  i,j\right)  \in\left\{
1,2,\ldots,n\right\}  ^{2};\\i<j\text{ and }\pi\left(  i\right)  <\pi\left(
j\right)  }}$\textquotedblright.]
\end{exercise}

Exercise \ref{exe.perm.pij-pii-sum} is \cite[Proposition 2.4]{SacUlf11}.

\begin{exercise}
\label{exe.perm.rearrangement}Whenever $m$ is an integer, we shall use the
notation $\left[  m\right]  $ for the set $\left\{  1,2,\ldots,m\right\}  $.
Also, recall Definition \ref{def.transpos.ii}.

Let $n\in\mathbb{N}$. Let $\sigma\in S_{n}$. Exercise \ref{exe.transpos.code}
shows that there is a unique $n$-tuple $\left(  i_{1},i_{2},\ldots
,i_{n}\right)  \in\left[  1\right]  \times\left[  2\right]  \times\cdots
\times\left[  n\right]  $ such that%
\[
\sigma=t_{1,i_{1}}\circ t_{2,i_{2}}\circ\cdots\circ t_{n,i_{n}}.
\]
Consider this $\left(  i_{1},i_{2},\ldots,i_{n}\right)  $.

For each $k\in\left\{  0,1,\ldots,n\right\}  $, we define a permutation
$\sigma_{k}\in S_{n}$ by $\sigma_{k}=t_{1,i_{1}}\circ t_{2,i_{2}}\circ
\cdots\circ t_{k,i_{k}}$.

For each $k\in\left[  n\right]  $, we let $m_{k}=\sigma_{k}\left(  k\right)  $.

\textbf{(a)} Prove that $\sigma_{k}\left(  i\right)  =i$ for each $i\in\left[
n\right]  $ and each $k\in\left\{  0,1,\ldots,i-1\right\}  $.

\textbf{(b)} Prove that $m_{k}\in\left[  k\right]  $ for all $k\in\left[
n\right]  $.

\textbf{(c)} Prove that $\sigma_{k}\left(  i_{k}\right)  =k$ for all
$k\in\left[  n\right]  $.

\textbf{(d)} Prove that $\sigma_{k}=t_{k,m_{k}}\circ\sigma_{k-1}$ for all
$k\in\left[  n\right]  $.

\textbf{(e)} Show that $\sigma^{-1}=t_{1,m_{1}}\circ t_{2,m_{2}}\circ
\cdots\circ t_{n,m_{n}}$.

\textbf{(f)} Let $x_{1},x_{2},\ldots,x_{n},y_{1},y_{2},\ldots,y_{n}$ be any
$2n$ numbers (e.g., rational numbers or real numbers or complex numbers).
Prove that
\[
\sum_{k=1}^{n}x_{k}y_{k}-\sum_{k=1}^{n}x_{k}y_{\sigma\left(  k\right)  }%
=\sum_{k=1}^{n}\left(  x_{i_{k}}-x_{k}\right)  \left(  y_{m_{k}}-y_{k}\right)
.
\]

\textbf{(g)} Now assume that the numbers $x_{1},x_{2},\ldots,x_{n},y_{1}%
,y_{2},\ldots,y_{n}$ are real and satisfy $x_{1}\geq x_{2}\geq\cdots\geq
x_{n}$ and $y_{1}\geq y_{2}\geq\cdots\geq y_{n}$. Conclude that
\[
\sum_{k=1}^{n}x_{k}y_{k}\geq\sum_{k=1}^{n}x_{k}y_{\sigma\left(  k\right)  }.
\]

\end{exercise}

\begin{remark}
The claim of Exercise \ref{exe.perm.rearrangement} \textbf{(g)} is known as
the \textit{rearrangement inequality}. It has several simple proofs (see,
e.g., \href{https://en.wikipedia.org/wiki/Rearrangement_inequality}{its
Wikipedia page}); the approach suggested by Exercise
\ref{exe.perm.rearrangement} is probably the most complicated, but it has the
advantage of giving an \textquotedblleft explicit\textquotedblright\ formula
for the difference between the two sides (in Exercise
\ref{exe.perm.rearrangement} \textbf{(f)}).
\end{remark}

\begin{exercise}
\label{exe.perm.order}Let $n\in\mathbb{N}$. Let $d=\operatorname{lcm}\left(
1,2,\ldots,n\right)  $. (Here, \textquotedblleft$\operatorname{lcm}%
$\textquotedblright\ stands for the
\href{https://en.wikipedia.org/wiki/Least_common_multiple}{least common
multiple} of several integers: Thus, $\operatorname{lcm}\left(  1,2,\ldots
,n\right)  $ is the smallest positive integer that is divisible by
$1,2,\ldots,n$.)

\textbf{(a)} Show that $\pi^{d}=\operatorname*{id}$ for every $\pi\in S_{n}$.

\textbf{(b)} Let $k$ be an integer such that every $\pi\in S_{n}$ satisfies
$\pi^{k}=\operatorname*{id}$. Show that $d\mid k$.
\end{exercise}

\begin{exercise}
\label{exe.perm.sigmacrosstau}Let $U$ and $V$ be two finite sets. Let $\sigma$
be a permutation of $U$. Let $\tau$ be a permutation of $V$. We define a
permutation $\sigma\times\tau$ of the set $U\times V$ by setting%
\[
\left(  \sigma\times\tau\right)  \left(  a,b\right)  =\left(  \sigma\left(
a\right)  ,\tau\left(  b\right)  \right)  \ \ \ \ \ \ \ \ \ \ \text{for every
}\left(  a,b\right)  \in U\times V.
\]

\textbf{(a)} Prove that $\sigma\times\tau$ is a well-defined permutation.

\textbf{(b)} Prove that $\sigma\times\tau=\left(  \sigma\times
\operatorname*{id}\nolimits_{V}\right)  \circ\left(  \operatorname*{id}%
\nolimits_{U}\times\tau\right)  $.

\textbf{(c)} Prove that $\left(  -1\right)  ^{\sigma\times\tau}=\left(
\left(  -1\right)  ^{\sigma}\right)  ^{\left\vert V\right\vert }\left(
\left(  -1\right)  ^{\tau}\right)  ^{\left\vert U\right\vert }$. (See Exercise
\ref{exe.ps4.2} for the definition of the signs $\left(  -1\right)
^{\sigma\times\tau}$, $\left(  -1\right)  ^{\sigma}$ and $\left(  -1\right)
^{\tau}$ appearing here.)
\end{exercise}

\begin{exercise}
\label{exe.perm.footrule}Let $n\in\mathbb{N}$. Let $\left[  n\right]  $ denote
the set $\left\{  1,2,\ldots,n\right\}  $. For each $\sigma\in S_{n}$, define
an integer $h\left(  \sigma\right)  $ by%
\[
h\left(  \sigma\right)  =\sum_{i\in\left[  n\right]  }\left\vert \sigma\left(
i\right)  -i\right\vert .
\]

Let $\sigma\in S_{n}$.

\textbf{(a)} Prove that%
\[
h\left(  \sigma\right)  =2\sum_{\substack{i\in\left[  n\right]  ;\\\sigma
\left(  i\right)  >i}}\left(  \sigma\left(  i\right)  -i\right)
=2\sum_{\substack{i\in\left[  n\right]  ;\\\sigma\left(  i\right)  <i}}\left(
i-\sigma\left(  i\right)  \right)  .
\]

\textbf{(b)} Prove that $h\left(  \sigma\circ\tau\right)  \leq h\left(
\sigma\right)  +h\left(  \tau\right)  $ for any $\tau\in S_{n}$.

\textbf{(c)} Prove that $h\left(  s_{k}\circ\sigma\right)  \leq h\left(
\sigma\right)  +2$ for each $k\in\left\{  1,2,\ldots,n-1\right\}  $.

\textbf{(d)} Prove that
\[
h\left(  \sigma\right)  /2\leq\ell\left(  \sigma\right)  \leq h\left(
\sigma\right)  .
\]

[\textbf{Hint:} The second inequality in part \textbf{(d)} is tricky. One way
to proceed is by classifying all inversions $\left(  i,j\right)  $ of $\sigma$
into two types: \textit{Type-I inversions} are those that satisfy
$\sigma\left(  i\right)  <j$, whereas \textit{Type-II inversions} are those
that satisfy $\sigma\left(  i\right)  \geq j$. Prove that the number of Type-I
inversions is $\leq\sum_{\substack{j\in\left[  n\right]  ;\\\sigma\left(
j\right)  <j}}\left(  j-\sigma\left(  j\right)  -1\right)  \leq\sum
_{\substack{i\in\left[  n\right]  ;\\\sigma\left(  i\right)  <i}}\left(
i-\sigma\left(  i\right)  \right)  $, whereas the number of Type-II inversions
is $\leq\sum_{\substack{i\in\left[  n\right]  ;\\\sigma\left(  i\right)
>i}}\left(  \sigma\left(  i\right)  -i\right)  $. Add these together to obtain
an upper bound on $\ell\left(  \sigma\right)  $.]
\end{exercise}

Exercise \ref{exe.perm.footrule} is a result of Diaconis and Graham
\cite[(3.5)]{DiaGra77}\footnote{Note that their notations are different; what
they call $I\left(  \pi\right)  $ and $D\left(  \pi\right)  $ would be called
$\ell\left(  \pi\right)  $ and $h\left(  \pi\right)  $ (respectively) in my
terminology.}. The integer $h\left(  \sigma\right)  $ defined in Exercise
\ref{exe.perm.footrule} is called \textit{Spearman's disarray} or
\textit{total displacement} of $\sigma$. A related concept (the depth of a
permutation) has been studied by Petersen and Tenner \cite{PetTen14}.

The next two exercises concern the inversions of a permutation. They use the
following definition:

\begin{definition}
\label{def.Inv}Let $n\in\mathbb{N}$. For every $\sigma\in S_{n}$, we let
$\operatorname*{Inv}\sigma$ denote the set of all inversions of $\sigma$.
\end{definition}

Exercise \ref{exe.ps2.2.5} \textbf{(c)} shows that any $n\in\mathbb{N}$ and
any two permutations $\sigma$ and $\tau$ in $S_{n}$ satisfy the inequality
$\ell\left(  \sigma\circ\tau\right)  \leq\ell\left(  \sigma\right)
+\ell\left(  \tau\right)  $. In the following exercise, we will see when this
inequality becomes an equality:

\begin{exercise}
\label{exe.perm.Inv.sub}Let $n\in\mathbb{N}$. Let $\sigma\in S_{n}$ and
$\tau\in S_{n}$.

\textbf{(a)} Prove that $\ell\left(  \sigma\circ\tau\right)  =\ell\left(
\sigma\right)  +\ell\left(  \tau\right)  $ holds if and only if
$\operatorname*{Inv}\tau\subseteq\operatorname*{Inv}\left(  \sigma\circ
\tau\right)  $.

\textbf{(b)} Prove that $\ell\left(  \sigma\circ\tau\right)  =\ell\left(
\sigma\right)  +\ell\left(  \tau\right)  $ holds if and only if
$\operatorname*{Inv}\left(  \sigma^{-1}\right)  \subseteq\operatorname*{Inv}%
\left(  \tau^{-1}\circ\sigma^{-1}\right)  $.

\textbf{(c)} Prove that $\operatorname*{Inv}\sigma\subseteq\operatorname*{Inv}%
\tau$ holds if and only if $\ell\left(  \tau\right)  =\ell\left(  \tau
\circ\sigma^{-1}\right)  +\ell\left(  \sigma\right)  $.

\textbf{(d)} Prove that if $\operatorname*{Inv}\sigma=\operatorname*{Inv}\tau
$, then $\sigma=\tau$.

\textbf{(e)} Prove that $\ell\left(  \sigma\circ\tau\right)  =\ell\left(
\sigma\right)  +\ell\left(  \tau\right)  $ holds if and only if $\left(
\operatorname*{Inv}\sigma\right)  \cap\left(  \operatorname*{Inv}\left(
\tau^{-1}\right)  \right)  =\varnothing$.
\end{exercise}

Exercise \ref{exe.perm.Inv.sub} \textbf{(d)} shows that if two permutations in
$S_{n}$ have the same set of inversions, then they are equal. In other words,
a permutation in $S_{n}$ is uniquely determined by its set of inversions. The
next exercise shows what set of inversions a permutation can have:

\begin{exercise}
\label{exe.perm.invset}Let $n\in\mathbb{N}$. Let $G=\left\{  \left(
i,j\right)  \in\mathbb{Z}^{2}\ \mid\ 1\leq i<j\leq n\right\}  $.

A subset $U$ of $G$ is said to be \textit{transitive} if every $a,b,c\in
\left\{  1,2,\ldots,n\right\}  $ satisfying $\left(  a,b\right)  \in U$ and
$\left(  b,c\right)  \in U$ also satisfy $\left(  a,c\right)  \in U$.

A subset $U$ of $G$ is said to be \textit{inversive} if there exists a
$\sigma\in S_{n}$ such that $U=\operatorname*{Inv}\sigma$.

Let $U$ be a subset of $G$. Prove that $U$ is inversive if and only if both
$U$ and $G\setminus U$ are transitive.
\end{exercise}

\section{\label{chp.det}An introduction to determinants}

In this chapter, we will define and study determinants in a combinatorial way
(in the spirit of Hefferon's book \cite{Hefferon}, Gill Williamson's notes
\cite[Chapter 3]{Gill}, Laue's notes \cite{Laue-det} and Zeilberger's paper
\cite{Zeilbe}\footnote{My notes differ from these sources in the following:
\par
\begin{itemize}
\item Hefferon's book \cite{Hefferon} is an introductory textbook for a first
course in Linear Algebra, and so treats rather little of the theory of
determinants (far less than what we do). It is, however, a good introduction
into the \textquotedblleft other part\textquotedblright\ of linear algebra
(i.e., the theory of vector spaces and linear maps), and puts determinants
into the context of that other part, which makes some of their properties
appear less mysterious. (Like many introductory textbooks, it only discusses
matrices over fields, not over commutative rings; it also uses more handwaving
in the proofs.)
\par
\item Zeilberger's paper \cite{Zeilbe} mostly proves advanced results (apart
from its Section 5, which proves our Theorem \ref{thm.det(AB)}). I would
recommend reading it after reading this chapter.
\par
\item Laue's notes \cite{Laue-det} are a brief introduction to determinants
that prove the main results in just 14 pages (although at the cost of terser
writing and stronger assumptions on the reader's preknowledge). If you read
these notes, make sure to pay attention to the \textquotedblleft Prerequisites
and some Terminology\textquotedblright\ section, as it explains the (unusual)
notations used in these notes.
\par
\item Gill Williamson's \cite[Chapter 3]{Gill} probably comes the closest to
what I am doing below (and is highly recommended, not least because it goes
much further into various interesting directions!). My notes are more
elementary and more detailed in what they do.
\end{itemize}
\par
Other references treating determinants in a combinatorial way are Day's
\cite[Chapter 6]{Day-proofs}, Herstein's \cite[\S 6.9]{Herstein}, Strickland's
\cite[\S 12 and Appendix B]{Strick13}, Mate's \cite{Mate14}, Walker's
\cite[\S 5.4]{Walker87}, and Pinkham's \cite[Chapter 11]{Pinkha15} (but they
all limit themselves to the basics).}). Nowadays, students usually learn about
determinants in the context of linear algebra, after having made the
acquaintance of vector spaces, matrices, linear transformations, Gaussian
elimination etc.; this approach to determinants (which I like to call the
\textquotedblleft linear-algebraic approach\textquotedblright) has certain
advantages and certain disadvantages compared to our combinatorial
approach\footnote{Its main advantage is that it gives more motivation and
context. However, the other (combinatorial) approach requires less
preknowledge and involves fewer technical subtleties (for example, it defines
the determinant directly by an explicit formula, while the linear-algebraic
approach defines it implicitly by a list of conditions which happen to
determine it uniquely), which is why I have chosen it.
\par
Examples of texts that introduce determinants via the linear-algebraic
approach are \cite[Chapter IX]{BirkMac}, \cite[Chapter 7]{GalQua18},
\cite[\S 8.3]{Goodman}, \cite[Chapter 5]{HoffmanKunze} and \cite[\S VII.3]%
{Hungerford-03}.
\par
Artin, in \cite[Chapter 1]{Artin}, takes a particularly quick approach to
determinants over a field (although it is quick at the cost of generality: for
example, the proof he gives for \cite[Theorem 1.4.9]{Artin} does not
generalize to matrices over commutative rings). Axler's \cite[Chapter
10B]{Axler} gives a singularly horrible treatment of determinants -- defining
them only for real and complex matrices and in a way that utterly hides their
combinatorial structure.}.

We shall study determinants of matrices over \textit{commutative
rings}.\footnote{This is a rather general setup, which includes determinants
of matrices with real entries, of matrices with complex entries, of matrices
with polynomial entries, and many other situations. One benefit of working
combinatorially is that studying determinants in this general setup is no more
difficult than studying them in more restricted settings.} First, let us
define what these words (\textquotedblleft commutative ring\textquotedblright,
\textquotedblleft matrix\textquotedblright\ and \textquotedblleft
determinant\textquotedblright) mean.

\subsection{\label{sect.commring}Commutative rings}

We begin by defining the concept of commutative rings, and exploring some
examples for it. This is only a brief introduction; much more can be found in
any text on abstract algebra (e.g., \cite{Artin}, \cite{Goodman},
\cite{Hungerford}, \cite{Hungerford-03}, \cite{Herstein} or \cite{ZarSam67}).

\begin{definition}
If $\mathbb{K}$ is a set, then a \textit{binary operation} on $\mathbb{K}$
means a map from $\mathbb{K}\times\mathbb{K}$ to $\mathbb{K}$. (In other
words, it means a function which takes two elements of $\mathbb{K}$ as input,
and returns an element of $\mathbb{K}$ as output.) For instance, the map from
$\mathbb{Z}\times\mathbb{Z}$ to $\mathbb{Z}$ which sends every pair $\left(
a,b\right)  \in\mathbb{Z}\times\mathbb{Z}$ to $3a-b$ is a binary operation on
$\mathbb{Z}$.

Sometimes, a binary operation $f$ on a set $\mathbb{K}$ will be
\textit{written infix}. This means that the image of $\left(  a,b\right)
\in\mathbb{K}\times\mathbb{K}$ under $f$ will be denoted by $afb$ instead of
$f\left(  a,b\right)  $. For instance, the binary operation $+$ on the set
$\mathbb{Z}$ (which sends a pair $\left(  a,b\right)  $ of integers to their
sum $a+b$) is commonly written infix, because one writes $a+b$ and not
$+\left(  a,b\right)  $ for the sum of $a$ and $b$.
\end{definition}

\begin{definition}
\label{def.commring}A \textit{commutative ring} means a set $\mathbb{K}$
endowed with

\begin{itemize}
\item two binary operations called \textquotedblleft
addition\textquotedblright\ and \textquotedblleft
multiplication\textquotedblright, and denoted by $+$ and $\cdot$,
respectively, and both written infix\footnotemark, and

\item two elements called $0_{\mathbb{K}}$ and $1_{\mathbb{K}}$
\end{itemize}

\noindent such that the following axioms are satisfied:

\begin{itemize}
\item \textit{Commutativity of addition:} We have $a+b=b+a$ for all
$a\in\mathbb{K}$ and $b\in\mathbb{K}$.

\item \textit{Commutativity of multiplication:} We have $ab=ba$ for all
$a\in\mathbb{K}$ and $b\in\mathbb{K}$. Here and in the following, the
expression \textquotedblleft$ab$\textquotedblright\ is shorthand for
\textquotedblleft$a\cdot b$\textquotedblright\ (as is usual for products of numbers).

\item \textit{Associativity of addition:} We have $a+\left(  b+c\right)
=\left(  a+b\right)  +c$ for all $a\in\mathbb{K}$, $b\in\mathbb{K}$ and
$c\in\mathbb{K}$.

\item \textit{Associativity of multiplication:} We have $a\left(  bc\right)
=\left(  ab\right)  c$ for all $a\in\mathbb{K}$, $b\in\mathbb{K}$ and
$c\in\mathbb{K}$.

\item \textit{Neutrality of }$0$\textit{:} We have $a+0_{\mathbb{K}%
}=0_{\mathbb{K}}+a=a$ for all $a\in\mathbb{K}$.

\item \textit{Existence of additive inverses:} For every $a\in\mathbb{K}$,
there exists an element $a^{\prime}\in\mathbb{K}$ such that $a+a^{\prime
}=a^{\prime}+a=0_{\mathbb{K}}$. This $a^{\prime}$ is commonly denoted by $-a$
and called the \textit{additive inverse} of $a$. (It is easy to check that it
is unique.)

\item \textit{Unitality (a.k.a. neutrality of }$1$\textit{):} We have
$1_{\mathbb{K}}a=a1_{\mathbb{K}}=a$ for all $a\in\mathbb{K}$.

\item \textit{Annihilation:} We have $0_{\mathbb{K}}a=a0_{\mathbb{K}%
}=0_{\mathbb{K}}$ for all $a\in\mathbb{K}$.

\item \textit{Distributivity:} We have $a\left(  b+c\right)  =ab+ac$ and
$\left(  a+b\right)  c=ac+bc$ for all $a\in\mathbb{K}$, $b\in\mathbb{K}$ and
$c\in\mathbb{K}$. Here and in the following, we are following the usual
convention (\textquotedblleft PEMDAS\textquotedblright) that
multiplication-like operations have higher precedence than addition-like
operations; thus, the expression \textquotedblleft$ab+ac$\textquotedblright%
\ must be understood as \textquotedblleft$\left(  ab\right)  +\left(
ac\right)  $\textquotedblright\ (and not, for example, as \textquotedblleft%
$a\left(  b+a\right)  c$\textquotedblright).
\end{itemize}
\end{definition}

\footnotetext{i.e., we write $a+b$ for the image of $\left(  a,b\right)
\in\mathbb{K}\times\mathbb{K}$ under the binary operation called
\textquotedblleft addition\textquotedblright, and we write $a\cdot b$ for the
image of $\left(  a,b\right)  \in\mathbb{K}\times\mathbb{K}$ under the binary
operation called \textquotedblleft multiplication\textquotedblright}(Some of
these axioms are redundant, in the sense that they can be derived from others.
For instance, the equality $\left(  a+b\right)  c=ac+bc$ can be derived from
the axiom $a\left(  b+c\right)  =ab+ac$ using commutativity of multiplication.
Also, annihilation follows from the other axioms\footnote{In fact, let
$a\in\mathbb{K}$. Distributivity yields $\left(  0_{\mathbb{K}}+0_{\mathbb{K}%
}\right)  a=0_{\mathbb{K}}a+0_{\mathbb{K}}a$, so that $0_{\mathbb{K}%
}a+0_{\mathbb{K}}a=\underbrace{\left(  0_{\mathbb{K}}+0_{\mathbb{K}}\right)
}_{\substack{=0_{\mathbb{K}}\\\text{(by neutrality of }0_{\mathbb{K}}\text{)}%
}}a=0_{\mathbb{K}}a$. Adding $-\left(  0_{\mathbb{K}}a\right)  $ on the left,
we obtain $-\left(  0_{\mathbb{K}}a\right)  +\left(  0_{\mathbb{K}%
}a+0_{\mathbb{K}}a\right)  =-\left(  0_{\mathbb{K}}a\right)  +0_{\mathbb{K}}%
a$. But $-\left(  0_{\mathbb{K}}a\right)  +0_{\mathbb{K}}a=0_{\mathbb{K}}$ (by
the definition of $-\left(  0_{\mathbb{K}}a\right)  $), and associativity of
addition shows that $-\left(  0_{\mathbb{K}}a\right)  +\left(  0_{\mathbb{K}%
}a+0_{\mathbb{K}}a\right)  =\underbrace{\left(  -\left(  0_{\mathbb{K}%
}a\right)  +0_{\mathbb{K}}a\right)  }_{=0_{\mathbb{K}}}+0_{\mathbb{K}%
}a=0_{\mathbb{K}}+0_{\mathbb{K}}a=0_{\mathbb{K}}a$ (by neutrality of
$0_{\mathbb{K}}$), so that $0_{\mathbb{K}}a=-\left(  0_{\mathbb{K}}a\right)
+\left(  0_{\mathbb{K}}a+0_{\mathbb{K}}a\right)  =-\left(  0_{\mathbb{K}%
}a\right)  +0_{\mathbb{K}}a=0_{\mathbb{K}}$. Thus, $0_{\mathbb{K}%
}a=0_{\mathbb{K}}$ is proven. Similarly one can show $a0_{\mathbb{K}%
}=0_{\mathbb{K}}$. Therefore, annihilation follows from the other axioms.}.
The reasons why we have chosen these axioms and not fewer (or more, or others)
are somewhat a matter of taste. For example, I like to explicitly require
annihilation, because it is an important axiom in the definition of a
\textit{semiring}, where it no longer follows from the others.)

\begin{definition}
As we have seen in Definition \ref{def.commring}, a commutative ring consists
of a set $\mathbb{K}$, two binary operations on this set named $+$ and $\cdot
$, and two elements of this set named $0$ and $1$. Thus, formally speaking, we
should encode a commutative ring as the $5$-tuple $\left(  \mathbb{K}%
,+,\cdot,0_{\mathbb{K}},1_{\mathbb{K}}\right)  $. Sometimes we will actually
do so; but most of the time, we will refer to the commutative ring just as the
\textquotedblleft commutative ring $\mathbb{K}$\textquotedblright, hoping that
the other four entries of the $5$-tuple (namely, $+$, $\cdot$, $0_{\mathbb{K}%
}$ and $1_{\mathbb{K}}$) are clear from the context. This kind of abbreviation
is commonplace in mathematics; it is called \textquotedblleft\textit{pars pro
toto}\textquotedblright\ (because we are referring to a large structure by the
same symbol as for a small part of it, and hoping that the rest can be
inferred from the context). It is an example of what is called
\textquotedblleft abuse of notation\textquotedblright.

The elements $0_{\mathbb{K}}$ and $1_{\mathbb{K}}$ of a commutative ring
$\mathbb{K}$ are called the \textit{zero} and the \textit{unity}%
\footnotemark\ of $\mathbb{K}$. They are usually denoted by $0$ and $1$
(without the subscript $\mathbb{K}$) when this can cause no confusion (and,
unfortunately, often also when it can). They are not always identical with the
actual integers $0$ and $1$.

The binary operations $+$ and $\cdot$ in Definition \ref{def.commring} are
also usually not identical with the binary operations $+$ and $\cdot$ on the
set of integers, and are denoted by $+_{\mathbb{K}}$ and $\cdot_{\mathbb{K}}$
when confusion can arise.

The set $\mathbb{K}$ is called the \textit{underlying set} of the commutative
ring $\mathbb{K}$. Let us again remind ourselves that the underlying set of a
commutative ring $\mathbb{K}$ is just a part of the data of $\mathbb{K}$.
\end{definition}

\footnotetext{Some people say \textquotedblleft unit\textquotedblright%
\ instead of \textquotedblleft unity\textquotedblright, but other people use
the word \textquotedblleft unit\textquotedblright\ for something different,
which makes every use of this word a potential pitfall.}Here are some examples
and non-examples of rings:\footnote{The following list of examples is long,
and some of these examples rely on knowledge that you might not have yet. As
usual with examples, you need not understand them all. When I say that Laurent
polynomial rings are examples of commutative rings, I do not assume that you
know what Laurent polynomials are; I merely want to ensure that, \textbf{if}
you have already encountered Laurent polynomials, then you get to know that
they form a commutative ring.}

\begin{itemize}
\item The sets $\mathbb{Z}$, $\mathbb{Q}$, $\mathbb{R}$ and $\mathbb{C}$
(endowed with the usual addition, the usual multiplication, the usual $0$ and
the usual $1$) are commutative rings. (Notice that existence of
\textbf{multiplicative} inverses is not required\footnote{A
\textit{multiplicative inverse} of an element $a\in\mathbb{K}$ means an
element $a^{\prime}\in\mathbb{K}$ such that $aa^{\prime}=a^{\prime
}a=1_{\mathbb{K}}$. (This is analogous to an additive inverse, except that
addition is replaced by multiplication, and $0_{\mathbb{K}}$ is replaced by
$1_{\mathbb{K}}$.) In a commutative ring, every element is required to have an
additive inverse (by the definition of a commutative ring), but not every
element is guaranteed to have a multiplicative inverse. (For instance, $2$ has
no multiplicative inverse in $\mathbb{Z}$, and $0$ has no multiplicative
inverse in any of the rings $\mathbb{Z}$, $\mathbb{Q}$, $\mathbb{R}$ and
$\mathbb{C}$.)
\par
We shall study multiplicative inverses more thoroughly in Section
\ref{sect.invertible} (where we will just call them \textquotedblleft
inverses\textquotedblright).}!)

\item The set $\mathbb{N}$ of nonnegative integers (again endowed with the
usual addition, the usual multiplication, the usual $0$ and the usual $1$) is
\textbf{not} a commutative ring. It fails the existence of additive inverses.
(Of course, negative numbers exist, but this does not count because they don't
lie in $\mathbb{N}$.)

\item We can define a commutative ring $\mathbb{Z}^{\prime}$ as follows: We
define a binary operation $\widetilde{\times}$ on $\mathbb{Z}$ (written infix)
by
\[
\left(  a\widetilde{\times}b=-ab\ \ \ \ \ \ \ \ \ \ \text{for all }\left(
a,b\right)  \in\mathbb{Z}\times\mathbb{Z}\right)  .
\]
Now, let $\mathbb{Z}^{\prime}$ be the \textbf{set} $\mathbb{Z}$, endowed with
the usual addition $+$ and the (unusual) multiplication $\widetilde{\times}$,
with the zero $0_{\mathbb{Z}^{\prime}}=0$ and with the unity $1_{\mathbb{Z}%
^{\prime}}=-1$. It is easy to check that $\mathbb{Z}^{\prime}$ is a
commutative ring\footnote{Notice that we have named this new commutative ring
$\mathbb{Z}^{\prime}$, not $\mathbb{Z}$ (despite having $\mathbb{Z}^{\prime
}=\mathbb{Z}$ as sets). The reason is that if we had named it $\mathbb{Z}$,
then we could no longer speak of \textquotedblleft the commutative ring
$\mathbb{Z}$\textquotedblright\ without being ambiguous (we would have to
specify every time whether we mean the usual multiplication or the unusual
one).}; it is an example of a commutative ring whose unity is clearly
\textbf{not} equal to the integer $1$ (which is why it is important to never
omit the subscript $\mathbb{Z}^{\prime}$ in $1_{\mathbb{Z}^{\prime}}$ here).

That said, $\mathbb{Z}^{\prime}$ is not a very interesting ring: It is
essentially \textquotedblleft a copy of $\mathbb{Z}$, except that every
integer $n$ has been renamed as $-n$\textquotedblright. To formalize this
intuition, we would need to introduce the notion of a \textit{ring
isomorphism}, which we don't want to do right here; but the idea is that the
bijection%
\[
\varphi:\mathbb{Z}\rightarrow\mathbb{Z}^{\prime},\ \ \ \ \ \ \ \ \ \ n\mapsto
-n
\]
satisfies%
\begin{align*}
\varphi\left(  a+b\right)   &  =\varphi\left(  a\right)  +\varphi\left(
b\right)  \ \ \ \ \ \ \ \ \ \ \text{for all }\left(  a,b\right)  \in
\mathbb{Z}\times\mathbb{Z};\\
\varphi\left(  a\cdot b\right)   &  =\varphi\left(  a\right)
\widetilde{\times}\varphi\left(  b\right)  \ \ \ \ \ \ \ \ \ \ \text{for all
}\left(  a,b\right)  \in\mathbb{Z}\times\mathbb{Z};\\
\varphi\left(  0\right)   &  =0_{\mathbb{Z}^{\prime}};\\
\varphi\left(  1\right)   &  =1_{\mathbb{Z}^{\prime}},
\end{align*}
and thus the ring $\mathbb{Z}^{\prime}$ can be viewed as the ring $\mathbb{Z}$
with its elements \textquotedblleft relabelled\textquotedblright\ using this bijection.

\item The polynomial rings $\mathbb{Z}\left[  x\right]  $, $\mathbb{Q}\left[
a,b\right]  $, $\mathbb{C}\left[  z_{1},z_{2},\ldots,z_{n}\right]  $ are
commutative rings. Laurent polynomial rings are also commutative rings. (Do
not worry if you have not seen these rings yet.)

\item The set of all functions $\mathbb{Q}\rightarrow\mathbb{Q}$ is a
commutative ring, where addition and multiplication are defined pointwise
(i.e., addition is defined by $\left(  f+g\right)  \left(  x\right)  =f\left(
x\right)  +g\left(  x\right)  $ for all $x\in\mathbb{Q}$, and multiplication
is defined by $\left(  fg\right)  \left(  x\right)  =f\left(  x\right)  \cdot
g\left(  x\right)  $ for all $x\in\mathbb{Q}$), where the zero is the
\textquotedblleft constant-$0$\textquotedblright\ function (sending every
$x\in\mathbb{Q}$ to $0$), and where the unity is the \textquotedblleft
constant-$1$\textquotedblright\ function (sending every $x\in\mathbb{Q}$ to
$1$). Of course, the same construction works if we consider functions
$\mathbb{R}\rightarrow\mathbb{C}$, or functions $\mathbb{C}\rightarrow
\mathbb{Q}$, or functions $\mathbb{N}\rightarrow\mathbb{Q}$, instead of
functions $\mathbb{Q}\rightarrow\mathbb{Q}$.\ \ \ \ \footnote{But not if we
consider functions $\mathbb{Q}\rightarrow\mathbb{N}$; such functions might
fail the existence of additive inverses.
\par
Generally, if $X$ is any set and $\mathbb{K}$ is any commutative ring, then
the set of all functions $X\rightarrow\mathbb{K}$ is a commutative ring, where
addition and multiplication are defined pointwise, where the zero is the
\textquotedblleft constant-$0_{\mathbb{K}}$\textquotedblright\ function, and
where the unity is the \textquotedblleft constant-$1_{\mathbb{K}}%
$\textquotedblright\ function.}

\item The set $\mathbb{S}$ of all real numbers of the form $a+b\sqrt{5}$ with
$a,b\in\mathbb{Q}$ (endowed with the usual notions of \textquotedblleft
addition\textquotedblright\ and \textquotedblleft
multiplication\textquotedblright\ defined on $\mathbb{R}$) is a commutative
ring\footnote{To prove this, we argue as follows:
\par
\begin{itemize}
\item Addition and multiplication are indeed two binary operations on
$\mathbb{S}$. This is because the sum of two elements of $\mathbb{S}$ is an
element of $\mathbb{S}$ (namely, $\left(  a+b\sqrt{5}\right)  +\left(
c+d\sqrt{5}\right)  =\left(  a+c\right)  +\left(  b+d\right)  \sqrt{5}$), and
so is their product (namely, $\left(  a+b\sqrt{5}\right)  \cdot\left(
c+d\sqrt{5}\right)  =\left(  ac+5bd\right)  +\left(  bc+ad\right)  \sqrt{5}$).
\par
\item All axioms of a commutative ring are satisfied for $\mathbb{S}$, except
maybe the existence of additive inverses. This is simply because the addition
and the multiplication in $\mathbb{S}$ are \textquotedblleft
inherited\textquotedblright\ from $\mathbb{R}$, and clearly all these axioms
come with the inheritance.
\par
\item Existence of additive inverses also holds in $\mathbb{S}$, because the
additive inverse of $a+b\sqrt{5}$ is $\left(  -a\right)  +\left(  -b\right)
\sqrt{5}$.
\end{itemize}
}.

\item We could define a different ring structure on the set $\mathbb{S}$ (that
is, a commutative ring which, as a set, is identical with $\mathbb{S}$, but
has a different choice of operations) as follows: We define a binary operation
$\ast$ on $\mathbb{S}$ by setting%
\[
\left(  a+b\sqrt{5}\right)  \ast\left(  c+d\sqrt{5}\right)  =ac+bd\sqrt
{5}\ \ \ \ \ \ \ \ \ \ \text{for all }\left(  a,b\right)  \in\mathbb{Q}%
\times\mathbb{Q}\text{ and }\left(  c,d\right)  \in\mathbb{Q}\times
\mathbb{Q}.
\]
\footnote{This is well-defined, because every element of $\mathbb{S}$ can be
written in the form $a+b\sqrt{5}$ for a \textbf{unique} pair $\left(
a,b\right)  \in\mathbb{Q}\times\mathbb{Q}$. This is a consequence of the
irrationality of $\sqrt{5}$.} Now, let $\mathbb{S}^{\prime}$ be the set
$\mathbb{S}$, endowed with the usual addition $+$ and the (unusual)
multiplication $\ast$, with the zero $0_{\mathbb{S}^{\prime}}=0$ and with the
unity $1_{\mathbb{S}^{\prime}}=1+\sqrt{5}$ (not the integer $1$). It is easy
to check that $\mathbb{S}^{\prime}$ is a commutative ring\footnote{Again, we
do not call it $\mathbb{S}$, in order to be able to distinguish between
different ring structures.}. The \textbf{sets} $\mathbb{S}$ and $\mathbb{S}%
^{\prime}$ are identical, but the \textbf{commutative rings} $\mathbb{S}$ and
$\mathbb{S}^{\prime}$ are not\footnote{Keep in mind that, due to our
\textquotedblleft pars pro toto\textquotedblright\ notation, \textquotedblleft
commutative ring $\mathbb{S}$\textquotedblright\ means more than
\textquotedblleft set $\mathbb{S}$\textquotedblright.}: For example, the ring
$\mathbb{S}^{\prime}$ has two nonzero elements whose product is $0$ (namely,
$1\ast\sqrt{5}=0$), whereas the ring $\mathbb{S}$ has no such things. This
shows that not only do we have $\mathbb{S}^{\prime}\neq\mathbb{S}$ as
commutative rings, but there is also no way to regard $\mathbb{S}^{\prime}$ as
\textquotedblleft a copy of $\mathbb{S}$ with its elements
renamed\textquotedblright\ (in the same way as we have regarded $\mathbb{Z}%
^{\prime}$ as \textquotedblleft a copy of $\mathbb{Z}$ with its elements
renamed\textquotedblright). This example should stress the point that a
commutative ring $\mathbb{K}$ is not just a set; it is a set endowed with two
operations ($+$ and $\cdot$) and two elements ($0_{\mathbb{K}}$ and
$1_{\mathbb{K}}$), and these operations and elements are no less important
than the set.

\item The set $\mathbb{S}_{3}$ of all real numbers of the form $a+b\sqrt[3]%
{5}$ with $a,b\in\mathbb{Q}$ (endowed with the usual addition, the usual
multiplication, the usual $0$ and the usual $1$) is \textbf{not} a commutative
ring. Indeed, multiplication is not a binary operation on this set
$\mathbb{S}_{3}$: It does not always send two elements of $\mathbb{S}_{3}$ to
an element of $\mathbb{S}_{3}$. For instance, $\left(  1+1\sqrt[3]{5}\right)
\left(  1+1\sqrt[3]{5}\right)  =1+2\sqrt[3]{5}+\left(  \sqrt[3]{5}\right)
^{2}$ is not in $\mathbb{S}_{3}$.

\item The set of all $2\times2$-matrices over $\mathbb{Q}$ is \textbf{not} a
commutative ring, because commutativity of multiplication does not hold for
this set. (In general, $AB\neq BA$ for matrices.)

\item If you like the empty set, you will enjoy the \textit{zero ring}. This
is the commutative ring which is defined as the one-element set $\left\{
0\right\}  $, with zero and unity both being $0$ (nobody said that they have
to be distinct!), with addition given by $0+0=0$ and with multiplication given
by $0\cdot0=0$. Of course, it is not an empty set\footnote{A commutative ring
cannot be empty, as it contains at least one element (namely, $0$).}, but it
plays a similar role in the world of commutative rings as the empty set does
in the world of sets: It carries no information itself, but things would break
if it were to be excluded\footnote{Some authors \textbf{do} prohibit the zero
ring from being a commutative ring (by requiring every commutative ring to
satisfy $0\neq1$). I think most of them run into difficulties from this
decision sooner or later.}.

Notice that the zero and the unity of the zero ring are identical, i.e., we
have $0_{\mathbb{K}}=1_{\mathbb{K}}$. This shows why it is dangerous to omit
the subscripts and just denote the zero and the unity by $0$ and $1$; in fact,
you don't want to rewrite the equality $0_{\mathbb{K}}=1_{\mathbb{K}}$ as
\textquotedblleft$0=1$\textquotedblright! (Most algebraists make a compromise
between wanting to omit the subscripts and having to clarify what $0$ and $1$
mean: They say that \textquotedblleft$0=1$ in $\mathbb{K}$\textquotedblright%
\ to mean \textquotedblleft$0_{\mathbb{K}}=1_{\mathbb{K}}$\textquotedblright.)

Generally, a \textit{trivial ring} is defined to be a commutative ring
containing only one element (which then necessarily is both the zero and the
unity of this ring). The addition and the multiplication of a trivial ring are
uniquely determined (since there is only one possible value that a sum or a
product could take). Every trivial ring can be viewed as the zero ring with
its element $0$ relabelled.\footnote{In more formal terms, the preceding
statement would say that \textquotedblleft every trivial ring is isomorphic to
the zero ring\textquotedblright.}

\item In set theory, the \textit{symmetric difference} of two sets $A$ and $B$
is defined to be the set $\left(  A\cup B\right)  \setminus\left(  A\cap
B\right)  =\left(  A\setminus B\right)  \cup\left(  B\setminus A\right)  $.
This symmetric difference is denoted by $A\bigtriangleup B$. Now, let $S$ be
any set. Let $\mathcal{P}\left(  S\right)  $ denote the powerset of $S$ (that
is, the set of all subsets of $S$). It is easy to check that the following ten
properties hold:%
\begin{align*}
A\bigtriangleup B  &  =B\bigtriangleup A\ \ \ \ \ \ \ \ \ \ \text{for any sets
}A\text{ and }B\text{;}\\
A\cap B  &  =B\cap A\ \ \ \ \ \ \ \ \ \ \text{for any sets }A\text{ and
}B\text{;}\\
\left(  A\bigtriangleup B\right)  \bigtriangleup C  &  =A\bigtriangleup\left(
B\bigtriangleup C\right)  \ \ \ \ \ \ \ \ \ \ \text{for any sets }A\text{,
}B\text{ and }C\text{;}\\
\left(  A\cap B\right)  \cap C  &  =A\cap\left(  B\cap C\right)
\ \ \ \ \ \ \ \ \ \ \text{for any sets }A\text{, }B\text{ and }C\text{;}\\
A\bigtriangleup\varnothing &  =\varnothing\bigtriangleup
A=A\ \ \ \ \ \ \ \ \ \ \text{for any set }A\text{;}\\
A\bigtriangleup A  &  =\varnothing\ \ \ \ \ \ \ \ \ \ \text{for any set
}A\text{;}\\
A\cap S  &  =S\cap A=A\ \ \ \ \ \ \ \ \ \ \text{for any subset }A\text{ of
}S\text{;}\\
\varnothing\cap A  &  =A\cap\varnothing=\varnothing
\ \ \ \ \ \ \ \ \ \ \text{for any set }A\text{;}\\
A\cap\left(  B\bigtriangleup C\right)   &  =\left(  A\cap B\right)
\bigtriangleup\left(  A\cap C\right)  \ \ \ \ \ \ \ \ \ \ \text{for any sets
}A\text{, }B\text{ and }C\text{;}\\
\left(  A\bigtriangleup B\right)  \cap C  &  =\left(  A\cap C\right)
\bigtriangleup\left(  B\cap C\right)  \ \ \ \ \ \ \ \ \ \ \text{for any sets
}A\text{, }B\text{ and }C\text{.}%
\end{align*}
Therefore, $\mathcal{P}\left(  S\right)  $ becomes a commutative ring, where
the addition is defined to be the operation $\bigtriangleup$, the
multiplication is defined to be the operation $\cap$, the zero is defined to
be the set $\varnothing$, and the unity is defined to be the set
$S$.\ \ \ \ \footnote{The ten properties listed above show that the axioms of
a commutative ring are satisfied for $\left(  \mathcal{P}\left(  S\right)
,\bigtriangleup,\cap,\varnothing,S\right)  $. In particular, the sixth
property shows that every subset $A$ of $S$ has an additive inverse -- namely,
itself. Of course, it is unusual for an element of a commutative ring to be
its own additive inverse, but in this example it happens all the time!}

The commutative ring $\mathcal{P}\left(  S\right)  $ has the property that
$a\cdot a=a$ for every $a\in\mathcal{P}\left(  S\right)  $. (This simply means
that $A\cap A=A$ for every $A\subseteq S$.) Commutative rings that have this
property are called
\href{https://en.wikipedia.org/wiki/Boolean_ring}{\textit{Boolean rings}}. (Of
course, $\mathcal{P}\left(  S\right)  $ is the eponymic example for a Boolean
ring; but there are also others.)

\item For every positive integer $n$, the residue classes of integers modulo
$n$ form a commutative ring, which is called $\mathbb{Z}/n\mathbb{Z}$ or
$\mathbb{Z}_{n}$ (depending on the author). This ring has $n$ elements (often
called \textquotedblleft integers modulo $n$\textquotedblright). When $n$ is a
composite number (e.g., $n=6$), this ring has the property that products of
nonzero\footnote{An element $a$ of a commutative ring $\mathbb{K}$ is said to
be \textit{nonzero} if $a\neq0_{\mathbb{K}}$. (This is not the same as saying
that $a$ is not the integer $0$, because the integer $0$ might not be
$0_{\mathbb{K}}$.)} elements can be zero (e.g., we have $2\cdot3\equiv
0\operatorname{mod}6$); this means that there is no way to define division by
all nonzero elements in this ring (even if we are allowed to create
fractions). Notice that $\mathbb{Z}/1\mathbb{Z}$ is a trivial ring.

We notice that if $n$ is a positive integer, and if $\mathbb{K}$ is the
commutative ring $\mathbb{Z}/n\mathbb{Z}$, then $\underbrace{1_{\mathbb{K}%
}+1_{\mathbb{K}}+\cdots+1_{\mathbb{K}}}_{n\text{ times}}=0_{\mathbb{K}}$
(because the left hand side of this equality is the residue class of $n$
modulo $n$, while the right hand side is the residue class of $0$ modulo $n$,
and these two residue classes are clearly equal).

\item Let us try to define \textquotedblleft division by
zero\textquotedblright. So, we introduce a new symbol $\infty$, and we try to
extend the addition on $\mathbb{Q}$ to the set $\mathbb{Q}\cup\left\{
\infty\right\}  $ by setting $a+\infty=\infty$ for all $a\in\mathbb{Q}%
\cup\left\{  \infty\right\}  $. We might also try to extend the multiplication
in some way, and perhaps to add some more elements (such as another symbol
$-\infty$ to serve as the product $\left(  -1\right)  \infty$). I claim that
(whatever we do with the multiplication, and whatever new elements we add) we
do not get a commutative ring. Indeed, assume the contrary. Thus, there exists
a commutative ring $\mathbb{W}$ which contains $\mathbb{Q}\cup\left\{
\infty\right\}  $ as a subset, and which has $a+\infty=\infty$ for all
$a\in\mathbb{Q}$. Thus, in $\mathbb{W}$, we have $1+\infty=\infty=0+\infty$.
Adding $\left(  -1\right)  \infty$ to both sides of this equality, we obtain
$1+\infty+\left(  -1\right)  \infty=0+\infty+\left(  -1\right)  \infty$, so
that $1=0$\ \ \ \ \footnote{because $\infty+\left(  -1\right)  \infty
=1\infty+\left(  -1\right)  \infty=\underbrace{\left(  1+\left(  -1\right)
\right)  }_{=0}\infty=0\infty=0$}; but this is absurd. Hence, we have found a
contradiction. This is why \textquotedblleft division by zero is
impossible\textquotedblright: One can define objects that behave like
\textquotedblleft infinity\textquotedblright\ (and they \textbf{are} useful),
but they break various standard rules such as the axioms of a commutative
ring. In contrast to this, adding a \textquotedblleft number\textquotedblright%
\ $i$ satisfying $i^{2}=-1$ to the real numbers is harmless: The complex
numbers $\mathbb{C}$ are still a commutative ring.

\item Here is an \textquotedblleft almost-ring\textquotedblright\ beloved to
many combinatorialists: the \textit{max-plus semiring} $\mathbb{T}$ (also
called the \textit{tropical semiring}\footnote{Caution: Both of these names
mean many other things as well.}). We create a new symbol $-\infty$, and we
set $\mathbb{T}=\mathbb{Z}\cup\left\{  -\infty\right\}  $ as sets, but we do
\textbf{not} \textquotedblleft inherit\textquotedblright\ the addition and the
multiplication from $\mathbb{Z}$. Instead, we denote the \textquotedblleft
addition\textquotedblright\ and \textquotedblleft
multiplication\textquotedblright\ operations on $\mathbb{Z}$ by $+_{\mathbb{Z}%
}$ and $\cdot_{\mathbb{Z}}$, and we define two new \textquotedblleft
addition\textquotedblright\ and \textquotedblleft
multiplication\textquotedblright\ operations $+_{\mathbb{T}}$ and
$\cdot_{\mathbb{T}}$ on $\mathbb{T}$ as follows:%
\begin{align*}
a+_{\mathbb{T}}b  &  =\max\left\{  a,b\right\}  ;\\
a\cdot_{\mathbb{T}}b  &  =a+_{\mathbb{Z}}b.
\end{align*}
(Here, we set $\max\left\{  -\infty,n\right\}  =\max\left\{  n,-\infty
\right\}  =n$ and $\left(  -\infty\right)  +_{\mathbb{Z}}n=n+_{\mathbb{Z}%
}\left(  -\infty\right)  =-\infty$ for every $n\in\mathbb{T}$.)

It turns out that the set $\mathbb{T}$ endowed with the two operations
$+_{\mathbb{T}}$ and $\cdot_{\mathbb{T}}$, the zero $0_{\mathbb{T}}=-\infty$
and the unity $1_{\mathbb{T}}=0$ comes rather close to being a commutative
ring. It satisfies all axioms of a commutative ring except for the existence
of additive inverses. Such a structure is called a \textit{semiring}. Other
examples of semirings are $\mathbb{N}$ and a reasonably defined $\mathbb{N}%
\cup\left\{  \infty\right\}  $ (with $0\infty=0$ and $a\infty=\infty$ for all
$a>0$).
\end{itemize}

If $\mathbb{K}$ is a commutative ring, then we can define a subtraction in
$\mathbb{K}$, even though we have not required a subtraction operation as part
of the definition of a commutative ring $\mathbb{K}$. Namely, the
\textit{subtraction} of a commutative ring $\mathbb{K}$ is the binary
operation $-$ on $\mathbb{K}$ (again written infix) defined as follows: For
every $a \in\mathbb{K}$ and $b \in\mathbb{K}$, set $a-b = a+b^{\prime}$, where
$b^{\prime}$ is the additive inverse of $b$. It is not hard to check that
$a-b$ is the unique element $c$ of $\mathbb{K}$ satisfying $a = b+c$; thus,
subtraction is ``the undoing of addition'' just as in the classical situation
of integers. Again, the notation $-$ for the subtraction of $\mathbb{K}$ is
denoted by $-_{\mathbb{K}}$ whenever a confusion with the subtraction of
integers could arise.

Whenever $a$ is an element of a commutative ring $\mathbb{K}$, we write $-a$
for the additive inverse of $a$. This is the same as $0_{\mathbb{K}} - a$.

The intuition for commutative rings is essentially that all computations that
can be performed with the operations $+$, $-$ and $\cdot$ on integers can be
similarly made in any commutative ring. For instance, if $a_{1},a_{2}%
,\ldots,a_{n}$ are $n$ elements of a commutative ring, then the sum
$a_{1}+a_{2}+\cdots+a_{n}$ is well-defined, and can be computed by adding the
elements $a_{1},a_{2},\ldots,a_{n}$ to each other in any order\footnote{For
instance, we can compute the sum $a+b+c+d$ of four elements $a,b,c,d$ in many
ways: For example, we can first add $a$ and $b$, then add $c$ and $d$, and
finally add the two results; alternatively, we can first add $a$ and $b$, then
add $d$ to the result, then add $c$ to the result. In a commutative ring, all
such ways lead to the same result.}. More generally: If $S$ is a finite set,
if $\mathbb{K}$ is a commutative ring, and if $\left(  a_{s}\right)  _{s\in
S}$ is a $\mathbb{K}$-valued $S$-family\footnote{See Definition
\ref{def.ind.families.fams} for the definition of this notion.}, then the sum
$\sum_{s\in S}a_{s}$ is defined in the same way as finite sums of numbers were
defined in Section \ref{sect.sums-repetitorium} (but with $\mathbb{A}$
replaced by $\mathbb{K}$, of course\footnote{and, consequently, $0$ replaced
by $0_{\mathbb{K}}$}); this definition is still legitimate\footnote{i.e., the
result does not depend on the choice of $t$ in (\ref{eq.sum.def.1})}, and
these finite sums of elements of $\mathbb{K}$ satisfy the same properties as
finite sums of numbers (see Section \ref{sect.sums-repetitorium} for these
properties). All this can be proven in the same way as it was proven for
numbers (in Section \ref{sect.ind.gen-com} and Section
\ref{sect.sums-repetitorium}). The same holds for finite products.
Furthermore, if $n$ is an integer and $a$ is an element of a commutative ring
$\mathbb{K}$, then we define an element $na$ of $\mathbb{K}$ by%
\[
na=%
\begin{cases}
\underbrace{a+a+\cdots+a}_{n\text{ addends}}, & \text{if }n\geq0;\\
-\left(  \underbrace{a+a+\cdots+a}_{-n\text{ addends}}\right)  , & \text{if
}n<0
\end{cases}
.
\]
\footnote{Notice that this definition of $na$ is \textbf{not} a particular
case of the product of two elements of $\mathbb{K}$, because $n$ is not an
element of $\mathbb{K}$.}

If $n$ is a nonnegative integer and $a$ is an element of a commutative ring
$\mathbb{K}$, then $a^{n}$ is a well-defined element of $\mathbb{K}$ (namely,
$a^{n}=\underbrace{a\cdot a\cdot\cdots\cdot a}_{n\text{ factors}}$). In
particular, applying this definition to $n=0$, we obtain%
\[
a^{0}=\underbrace{a\cdot a\cdot\cdots\cdot a}_{0\text{ factors}}=\left(
\text{empty product}\right)  =1\ \ \ \ \ \ \ \ \ \ \text{for each }%
a\in\mathbb{K}.
\]

The following identities hold:%
\begin{align}
\left(  n+m\right)  a  &  =na+ma\ \ \ \ \ \ \ \ \ \ \text{for }a\in
\mathbb{K}\text{ and }n,m\in\mathbb{Z};\label{eq.rings.nmdist1}\\
n\left(  a+b\right)   &  =na+nb\ \ \ \ \ \ \ \ \ \ \text{for }a,b\in
\mathbb{K}\text{ and }n\in\mathbb{Z};\label{eq.rings.nmdist2}\\
-\left(  a+b\right)   &  =\left(  -a\right)  +\left(  -b\right)
\ \ \ \ \ \ \ \ \ \ \text{for }a,b\in\mathbb{K};\label{eq.rings.-(a+b)}\\
1a  &  =a\ \ \ \ \ \ \ \ \ \ \text{for }a\in\mathbb{K};\label{eq.rings.1a}\\
0a  &  =0_{\mathbb{K}}\ \ \ \ \ \ \ \ \ \ \text{for }a\in\mathbb{K}%
\label{eq.rings.0a}\\
&  \ \ \ \ \ \ \ \ \ \ \left(  \text{here, the \textquotedblleft%
}0\text{\textquotedblright\ on the left hand side means the integer }0\right)
;\nonumber\\
\left(  -1\right)  a  &  =-a\ \ \ \ \ \ \ \ \ \ \text{for }a\in\mathbb{K}%
;\label{eq.rings.-1a}\\
-\left(  -a\right)   &  =a\ \ \ \ \ \ \ \ \ \ \text{for }a\in\mathbb{K}%
;\label{eq.rings.-(-a)}\\
-\left(  ab\right)   &  =\left(  -a\right)  b=a\left(  -b\right)
\ \ \ \ \ \ \ \ \ \ \text{for }a,b\in\mathbb{K};\label{eq.rings.-(ab)}\\
-\left(  na\right)   &  =\left(  -n\right)  a=n\left(  -a\right)
\ \ \ \ \ \ \ \ \ \ \text{for }a\in\mathbb{K}\text{ and }n\in\mathbb{Z}%
;\label{eq.rings.-(na)}\\
n\left(  ab\right)   &  =\left(  na\right)  b=a\left(  nb\right)
\ \ \ \ \ \ \ \ \ \ \text{for }a,b\in\mathbb{K}\text{ and }n\in\mathbb{Z}%
;\label{eq.rings.nab}\\
\left(  nm\right)  a  &  =n\left(  ma\right)  \ \ \ \ \ \ \ \ \ \ \text{for
}a\in\mathbb{K}\text{ and }n,m\in\mathbb{Z};\label{eq.rings.nma}\\
n0_{\mathbb{K}}  &  =0_{\mathbb{K}}\ \ \ \ \ \ \ \ \ \ \text{for }%
n\in\mathbb{Z};\nonumber\\
1^{n}  &  =1\ \ \ \ \ \ \ \ \ \ \text{for }n\in\mathbb{N};\nonumber\\
0^{n}  &  =%
\begin{cases}
0, & \text{if }n>0;\\
1, & \text{if }n=0
\end{cases}
\ \ \ \ \ \ \ \ \ \ \text{for }n\in\mathbb{N};\label{eq.rings.0**n}\\
a^{n+m}  &  =a^{n}a^{m}\ \ \ \ \ \ \ \ \ \ \text{for }a\in\mathbb{K}\text{ and
}n,m\in\mathbb{N};\label{eq.rings.a**(n+m)}\\
a^{nm}  &  =\left(  a^{n}\right)  ^{m}\ \ \ \ \ \ \ \ \ \ \text{for }%
a\in\mathbb{K}\text{ and }n,m\in\mathbb{N};\nonumber\\
\left(  ab\right)  ^{n}  &  =a^{n}b^{n}\ \ \ \ \ \ \ \ \ \ \text{for }%
a,b\in\mathbb{K}\text{ and }n\in\mathbb{N}. \label{eq.rings.-(ab)**n}%
\end{align}
Here, we are using the standard notations $+$, $\cdot$, $0$ and $1$ for the
addition, the multiplication, the zero and the unity of $\mathbb{K}$, because
confusion (e.g., confusion of the $0$ with the integer $0$) is rather
unlikely.\footnote{For instance, in the statement \textquotedblleft$-\left(
a+b\right)  =\left(  -a\right)  +\left(  -b\right)  $ for $a,b\in\mathbb{K}%
$\textquotedblright, it is clear that the $+$ can only stand for the addition
of $\mathbb{K}$ and not (say) for the addition of integers (since $a$, $b$,
$-a$ and $-b$ are elements of $\mathbb{K}$, not (generally) integers). The
only statement whose meaning is ambiguous is \textquotedblleft$0^{n}=%
\begin{cases}
0, & \text{if }n>0;\\
1, & \text{if }n=0
\end{cases}
$ for $n\in\mathbb{N}$\textquotedblright. In this statement, the
\textquotedblleft$0$\textquotedblright\ in \textquotedblleft$n>0$%
\textquotedblright\ and the \textquotedblleft$0$\textquotedblright\ in
\textquotedblleft$n=0$\textquotedblright\ clearly mean the integer $0$ (since
they are being compared with the integer $n$), but the other two appearances
of \textquotedblleft$0$\textquotedblright\ and the \textquotedblleft%
$1$\textquotedblright\ are ambiguous. I hope that the context makes it clear
enough that they mean the zero and the unity of $\mathbb{K}$ (and not the
integers $0$ and $1$), because otherwise this equality would not be a
statement about $\mathbb{K}$ at all.} We shall keep doing so in the following,
apart from situations where confusion can realistically occur.\footnote{Notice
that the equalities (\ref{eq.rings.nab}) and (\ref{eq.rings.nma}) are
\textbf{not} particular cases of the associativity of multiplication which we
required to hold for $\mathbb{K}$. Indeed, the latter associativity says that
$a\left(  bc\right)  =\left(  ab\right)  c$ for all $a\in\mathbb{K}$,
$b\in\mathbb{K}$ and $c\in\mathbb{K}$. But in (\ref{eq.rings.nab}) and
(\ref{eq.rings.nma}), the $n$ is an integer, not an element of $\mathbb{K}$.}

The identities listed above are not hard to prove. Indeed, they are
generalizations of well-known identities holding for rational numbers; and
some of them (for example, (\ref{eq.rings.a**(n+m)}) and
(\ref{eq.rings.-(ab)**n})) can be proved in exactly the same way as those
identities for rational numbers.\footnote{For example, it is well-known that
\[
\left(  ab\right)  ^{n}=a^{n}b^{n}\ \ \ \ \ \ \ \ \ \ \text{for any }%
a,b\in\mathbb{Q}\text{ and }n\in\mathbb{N}.
\]
This can be easily proven by induction on $n$, using the commutativity and
associativity rules for multiplication of rational numbers and the fact that
$1\cdot1=1$. The same argument can be used to prove (\ref{eq.rings.-(ab)**n}).
The only change required is replacing every appearance of \textquotedblleft%
$\mathbb{Q}$\textquotedblright\ by \textquotedblleft$\mathbb{K}$%
\textquotedblright.}

If $a$ and $b$ are two elements of a commutative ring $\mathbb{K}$, then the
expression \textquotedblleft$-ab$\textquotedblright\ appears ambiguous, since
it can be interpreted either as \textquotedblleft$-\left(  ab\right)
$\textquotedblright\ or as \textquotedblleft$\left(  -a\right)  b$%
\textquotedblright. But (\ref{eq.rings.-(ab)}) shows that these two
interpretations yield the same result; thus, we can write this expression
\textquotedblleft$-ab$\textquotedblright\ without fearing ambiguity.
Similarly, if $n\in\mathbb{Z}$ and $a,b\in\mathbb{K}$, then the expression
\textquotedblleft$nab$\textquotedblright\ is unambiguous, because
(\ref{eq.rings.nab}) shows that the two possible ways to interpret it (namely,
as \textquotedblleft$n\left(  ab\right)  $\textquotedblright\ and as
\textquotedblleft$\left(  na\right)  b$\textquotedblright) yield the same
result. Similarly, if $n,m\in\mathbb{Z}$ and $a\in\mathbb{K}$, then the
expression \textquotedblleft$nma$\textquotedblright\ is unambiguous, because
of (\ref{eq.rings.nma}).

Furthermore, finite sums such as $\sum_{s\in S}a_{s}$ (where $S$ is a finite
set, and $a_{s}\in\mathbb{K}$ for every $s\in S$), and finite products such as
$\prod_{s\in S}a_{s}$ (where $S$ is a finite set, and $a_{s}\in\mathbb{K}$ for
every $s\in S$) are defined whenever $\mathbb{K}$ is a commutative ring.
Again, the definition is the same as for numbers, and these sums and products
behave as they do for numbers.\footnote{Of course, empty sums of elements of
$\mathbb{K}$ are defined to equal $0_{\mathbb{K}}$, and empty products of
elements of $\mathbb{K}$ are defined to equal $1_{\mathbb{K}}$.} For example,
Exercise \ref{exe.perm.sign.pseudoexplicit} still holds if we replace
\textquotedblleft$\mathbb{C}$\textquotedblright\ by \textquotedblleft%
$\mathbb{K}$\textquotedblright\ in it (and the same solution proves it)
whenever $\mathbb{K}$ is a commutative ring. From the fact that finite sums
and finite products of elements of $\mathbb{K}$ are well-defined, we can also
conclude that expressions such as \textquotedblleft$a_{1}+a_{2}+\cdots+a_{k}%
$\textquotedblright\ and \textquotedblleft$a_{1}a_{2}\cdots a_{k}%
$\textquotedblright\ (where $a_{1},a_{2},\ldots,a_{k}$ are finitely many
elements of $\mathbb{K}$) are well-defined.

Various identities that hold for numbers also hold for elements of arbitrary
commutative rings. For example, an analogue of the binomial formula
(Proposition \ref{prop.binom.binomial}) holds: If $\mathbb{K}$ is a
commutative ring, then%
\begin{equation}
\left(  a+b\right)  ^{n}=\sum_{k=0}^{n}\dbinom{n}{k}a^{k}b^{n-k}%
\ \ \ \ \ \ \ \ \ \ \text{for }a,b\in\mathbb{K}\text{ and }n\in\mathbb{N}.
\label{eq.rings.(a+b)**n}%
\end{equation}
(We can obtain a proof of (\ref{eq.rings.(a+b)**n}) by re-reading the solution
to Exercise \ref{exe.prop.binom.binomial}, while replacing every
\textquotedblleft$x$\textquotedblright\ by an \textquotedblleft$a$%
\textquotedblright\ and replacing every \textquotedblleft$y$\textquotedblright%
\ by a \textquotedblleft$b$\textquotedblright. Another proof of
(\ref{eq.rings.(a+b)**n}) is given in the solution to Exercise
\ref{exe.prod(ai+bi)} \textbf{(b)}.)

\begin{remark}
The notion of a \textquotedblleft commutative ring\textquotedblright\ is not
fully standardized; there exist several competing definitions:

For some people, a \textquotedblleft commutative ring\textquotedblright\ is
\textbf{not} endowed with an element $1$ (although it \textbf{can} have such
an element), and, consequently, does not have to satisfy the unitality axiom.
According to their definition, for example, the set
\[
\left\{  \ldots,-4,-2,0,2,4,\ldots\right\}  =\left\{  2n\ \mid\ n\in
\mathbb{Z}\right\}  =\left(  \text{the set of all even integers}\right)
\]
is a commutative ring (with the usual addition and multiplication). (In
contrast, our definition of a \textquotedblleft commutative
ring\textquotedblright\ does not accept this set as a commutative ring,
because it does not contain any element which would fill the role of $1$.)
These people tend to use the notation \textquotedblleft commutative ring with
unity\textquotedblright\ (or \textquotedblleft commutative ring with
$1$\textquotedblright) to mean a commutative ring which is endowed with a $1$
and satisfies the unitality axiom (i.e., what we call a \textquotedblleft
commutative ring\textquotedblright).

On the other hand, there are authors who use the word \textquotedblleft
ring\textquotedblright\ for what we call \textquotedblleft commutative
ring\textquotedblright. These are mostly the authors who work with commutative
rings all the time and find the name \textquotedblleft commutative
ring\textquotedblright\ too long.

When you are reading about rings, it is important to know which meaning of
\textquotedblleft ring\textquotedblright\ the author is subscribing to. (Often
this can be inferred from the examples given.)
\end{remark}

\begin{exercise}
\label{exe.prod(ai+bi)}Let $\mathbb{K}$ be a commutative ring. For every
$n\in\mathbb{N}$, let $\left[  n\right]  $ denote the set $\left\{
1,2,\ldots,n\right\}  $.

\textbf{(a)} Let $n\in\mathbb{N}$. Let $a_{1},a_{2},\ldots,a_{n}$ be $n$
elements of $\mathbb{K}$. Let $b_{1},b_{2},\ldots,b_{n}$ be $n$ further
elements of $\mathbb{K}$. Prove that%
\[
\prod_{i=1}^{n}\left(  a_{i}+b_{i}\right)  =\sum_{I\subseteq\left[  n\right]
}\left(  \prod_{i\in I}a_{i}\right)  \left(  \prod_{i\in\left[  n\right]
\setminus I}b_{i}\right)  .
\]
(Here, as usual, the summation sign $\sum_{I\subseteq\left[  n\right]  }$
means $\sum_{I\in\mathcal{P}\left(  \left[  n\right]  \right)  }$, where
$\mathcal{P}\left(  \left[  n\right]  \right)  $ denotes the powerset of
$\left[  n\right]  $.)

\textbf{(b)} Use Exercise \ref{exe.prod(ai+bi)} to give a new proof of
(\ref{eq.rings.(a+b)**n}).
\end{exercise}

\begin{exercise}
\label{exe.multinom2}For each $m\in\mathbb{N}$ and $\left(  k_{1},k_{2}%
,\ldots,k_{m}\right)  \in\mathbb{N}^{m}$, let us define a positive integer
$\mathbf{m}\left(  k_{1},k_{2},\ldots,k_{m}\right)  $ by $\mathbf{m}\left(
k_{1},k_{2},\ldots,k_{m}\right)  =\dfrac{\left(  k_{1}+k_{2}+\cdots
+k_{m}\right)  !}{k_{1}!k_{2}!\cdots k_{m}!}$. (This is indeed a positive
integer, because Exercise \ref{exe.multinom1} says so.)

Let $\mathbb{K}$ be a commutative ring. Let $m\in\mathbb{N}$. Let $a_{1}%
,a_{2},\ldots,a_{m}$ be $m$ elements of $\mathbb{K}$. Let $n\in\mathbb{N}$.
Prove that%
\[
\left(  a_{1}+a_{2}+\cdots+a_{m}\right)  ^{n}=\sum_{\substack{\left(
k_{1},k_{2},\ldots,k_{m}\right)  \in\mathbb{N}^{m};\\k_{1}+k_{2}+\cdots
+k_{m}=n}}\mathbf{m}\left(  k_{1},k_{2},\ldots,k_{m}\right)  \prod_{i=1}%
^{m}a_{i}^{k_{i}}.
\]
(This is called the \textit{multinomial formula}.)
\end{exercise}

\subsection{Matrices}

We have briefly defined determinants in Definition \ref{def.det.old}, but we
haven't done much with them. This will be amended now. But let us first recall
the definitions of basic notions in matrix algebra.

In the following, we fix a commutative ring $\mathbb{K}$. The elements of
$\mathbb{K}$ will be called \textit{scalars} (to distinguish them from
\textit{vectors} and \textit{matrices}, which we will soon discuss, and which
are structures containing several elements of $\mathbb{K}$).

If you feel uncomfortable with commutative rings, you are free to think that
$\mathbb{K}=\mathbb{Q}$ or $\mathbb{K}=\mathbb{C}$ in the following; but
everything I am doing works for any commutative ring unless stated otherwise.

Given two nonnegative integers $n$ and $m$, an $n\times m$\textit{-matrix}
(or, more precisely, $n\times m$\textit{-matrix over} $\mathbb{K}$) means a
rectangular table with $n$ rows and $m$ columns whose entries are elements of
$\mathbb{K}$.\ \ \ \ \footnote{Formally speaking, this means that an $n\times
m$-matrix is a map from $\left\{  1,2,\ldots,n\right\}  \times\left\{
1,2,\ldots,m\right\}  $ to $\mathbb{K}$. We represent such a map as a
rectangular table by writing the image of $\left(  i,j\right)  \in\left\{
1,2,\ldots,n\right\}  \times\left\{  1,2,\ldots,m\right\}  $ into the cell in
the $i$-th row and the $j$-th column.
\par
Thus, the notion of an $n\times m$-matrix is closely akin to what we called an
\textquotedblleft$n\times m$-table of elements of $\mathbb{K}$%
\textquotedblright\ in Definition \ref{def.ind.families.rectab}. The main
difference between these two notions is that an $n\times m$-matrix
\textquotedblleft knows\textquotedblright\ $\mathbb{K}$, whereas an $n\times
m$-table does not (i.e., two $n\times m$-matrices that have the same entries
in the same positions but are defined using different commutative rings
$\mathbb{K}$ are considered different, but two such $n\times m$-tables are
considered identical).} For instance, when $\mathbb{K}=\mathbb{Q}$, the table
$\left(
\begin{array}
[c]{ccc}%
1 & -2/5 & 4\\
1/3 & -1/2 & 0
\end{array}
\right)  $ is a $2\times3$-matrix. A \textit{matrix} simply means an $n\times
m$-matrix for some $n\in\mathbb{N}$ and $m\in\mathbb{N}$. These $n$ and $m$
are said to be the \textit{dimensions} of the matrix.

If $A$ is an $n\times m$-matrix, and if $i\in\left\{  1,2,\ldots,n\right\}  $
and $j\in\left\{  1,2,\ldots,m\right\}  $, then the $\left(  i,j\right)
$\textit{-th entry of }$A$ means the entry of $A$ in row $i$ and column $j$.
For instance, the $\left(  1,2\right)  $-th entry of the matrix $\left(
\begin{array}
[c]{ccc}%
1 & -2/5 & 4\\
1/3 & -1/2 & 0
\end{array}
\right)  $ is $-2/5$.

If $n\in\mathbb{N}$ and $m\in\mathbb{N}$, and if we are given an element
$a_{i,j}\in\mathbb{K}$ for every $\left(  i,j\right)  \in\left\{
1,2,\ldots,n\right\}  \times\left\{  1,2,\ldots,m\right\}  $, then we use the
notation $\left(  a_{i,j}\right)  _{1\leq i\leq n,\ 1\leq j\leq m}$ for the
$n\times m$-matrix whose $\left(  i,j\right)  $-th entry is $a_{i,j}$ for all
$\left(  i,j\right)  \in\left\{  1,2,\ldots,n\right\}  \times\left\{
1,2,\ldots,m\right\}  $. Thus,%
\[
\left(  a_{i,j}\right)  _{1\leq i\leq n,\ 1\leq j\leq m}=\left(
\begin{array}
[c]{cccc}%
a_{1,1} & a_{1,2} & \cdots & a_{1,m}\\
a_{2,1} & a_{2,2} & \cdots & a_{2,m}\\
\vdots & \vdots & \ddots & \vdots\\
a_{n,1} & a_{n,2} & \cdots & a_{n,m}%
\end{array}
\right)  .
\]
The letters $i$ and $j$ are not set in stone; they are bound variables like
the $k$ in \textquotedblleft$\sum_{k=1}^{n}k$\textquotedblright. Thus, you are
free to write $\left(  a_{x,y}\right)  _{1\leq x\leq n,\ 1\leq y\leq m}$ or
$\left(  a_{j,i}\right)  _{1\leq j\leq n,\ 1\leq i\leq m}$ instead of $\left(
a_{i,j}\right)  _{1\leq i\leq n,\ 1\leq j\leq m}$ (and we will use this
freedom eventually).\footnote{Many authors love to abbreviate
\textquotedblleft$a_{i,j}$\textquotedblright\ by \textquotedblleft$a_{ij}%
$\textquotedblright\ (hoping that the reader will not mistake the subscript
\textquotedblleft$ij$\textquotedblright\ for a product or (in the case where
$i$ and $j$ are single-digit numbers) for a two-digit number). The only
advantage of this abbreviation that I am aware of is that it saves you a
comma; I do not understand why it is so popular. But you should be aware of it
in case you are reading other texts.}

Matrices can be added if they share the same dimensions: If $n$ and $m$ are
two nonnegative integers, and if $A=\left(  a_{i,j}\right)  _{1\leq i\leq
n,\ 1\leq j\leq m}$ and $B=\left(  b_{i,j}\right)  _{1\leq i\leq n,\ 1\leq
j\leq m}$ are two $n\times m$-matrices, then $A+B$ means the $n\times
m$-matrix $\left(  a_{i,j}+b_{i,j}\right)  _{1\leq i\leq n,\ 1\leq j\leq m}$.
Thus, matrices are added \textquotedblleft entry by entry\textquotedblright;
for example, $\left(
\begin{array}
[c]{ccc}%
a & b & c\\
d & e & f
\end{array}
\right)  +\left(
\begin{array}
[c]{ccc}%
a^{\prime} & b^{\prime} & c^{\prime}\\
d^{\prime} & e^{\prime} & f^{\prime}%
\end{array}
\right)  =\left(
\begin{array}
[c]{ccc}%
a+a^{\prime} & b+b^{\prime} & c+c^{\prime}\\
d+d^{\prime} & e+e^{\prime} & f+f^{\prime}%
\end{array}
\right)  $. Similarly, subtraction is defined: If $A=\left(  a_{i,j}\right)
_{1\leq i\leq n,\ 1\leq j\leq m}$ and $B=\left(  b_{i,j}\right)  _{1\leq i\leq
n,\ 1\leq j\leq m}$, then $A-B=\left(  a_{i,j}-b_{i,j}\right)  _{1\leq i\leq
n,\ 1\leq j\leq m}$.

Similarly, one can define the product of a scalar $\lambda\in\mathbb{K}$ with
a matrix $A$: If $\lambda\in\mathbb{K}$ is a scalar, and if $A=\left(
a_{i,j}\right)  _{1\leq i\leq n,\ 1\leq j\leq m}$ is an $n\times m$-matrix,
then $\lambda A$ means the $n\times m$-matrix $\left(  \lambda a_{i,j}\right)
_{1\leq i\leq n,\ 1\leq j\leq m}$.

Defining the product of two matrices is trickier. Matrices are \textbf{not}
multiplied \textquotedblleft entry by entry\textquotedblright; this would not
be a very interesting definition. Instead, their product is defined as
follows: If $n$, $m$ and $\ell$ are three nonnegative integers, then the
product $AB$ of an $n\times m$-matrix $A=\left(  a_{i,j}\right)  _{1\leq i\leq
n,\ 1\leq j\leq m}$ with an $m\times\ell$-matrix $B=\left(  b_{i,j}\right)
_{1\leq i\leq m,\ 1\leq j\leq\ell}$ means the $n\times\ell$-matrix%
\[
\left(  \sum_{k=1}^{m}a_{i,k}b_{k,j}\right)  _{1\leq i\leq n,\ 1\leq j\leq
\ell}.
\]
This definition looks somewhat counterintuitive, so let me comment on it.
First of all, for $AB$ to be defined, $A$ and $B$ are \textbf{not} required to
have the same dimensions; instead, $A$ must have as many columns as $B$ has
rows. The resulting matrix $AB$ then has as many rows as $A$ and as many
columns as $B$. Every entry of $AB$ is a sum of products of an entry of $A$
with an entry of $B$ (not a single such product). More precisely, the $\left(
i,j\right)  $-th entry of $AB$ is a sum of products of an entry in the $i$-th
row of $A$ with the respective entry in the $j$-th column of $B$. For example,%
\[
\left(
\begin{array}
[c]{ccc}%
a & b & c\\
d & e & f
\end{array}
\right)  \left(
\begin{array}
[c]{ccc}%
a^{\prime} & d^{\prime} & g^{\prime}\\
b^{\prime} & e^{\prime} & h^{\prime}\\
c^{\prime} & f^{\prime} & i^{\prime}%
\end{array}
\right)  =\left(
\begin{array}
[c]{ccc}%
aa^{\prime}+bb^{\prime}+cc^{\prime} & ad^{\prime}+be^{\prime}+cf^{\prime} &
ag^{\prime}+bh^{\prime}+ci^{\prime}\\
da^{\prime}+eb^{\prime}+fc^{\prime} & dd^{\prime}+ee^{\prime}+ff^{\prime} &
dg^{\prime}+eh^{\prime}+fi^{\prime}%
\end{array}
\right)  .
\]

The multiplication of matrices is not commutative! It is easy to find examples
of two matrices $A$ and $B$ for which the products $AB$ and $BA$ are distinct,
or one of them is well-defined but the other is not\footnote{This happens if
$A$ has as many columns as $B$ has rows, but $B$ does not have as many columns
as $A$ has rows.}.

For given $n\in\mathbb{N}$ and $m\in\mathbb{N}$, we define the $n\times
m$\textit{ zero matrix} to be the $n\times m$-matrix whose all entries are $0$
(that is, the $n\times m$-matrix $\left(  0\right)  _{1\leq i\leq n,\ 1\leq
j\leq m}$). We denote this matrix by $0_{n\times m}$. If $A$ is any $n\times
m$-matrix, then the $n\times m$-matrix $-A$ is defined to be $0_{n\times m}-A$.

A sum $\sum_{i\in I}A_{i}$ of finitely many matrices $A_{i}$ is defined in the
same way as a sum of numbers or of elements of a commutative
ring\footnote{with the caveat that an empty sum of $n\times m$-matrices is not
the number $0$, but the $n\times m$-matrix $0_{n,m}$}. However, a product
$\prod_{i\in I}A_{i}$ of finitely many matrices $A_{i}$ (in general) cannot be
defined, because the result would depend on the order of multiplication.

For every $n\in\mathbb{N}$, we let $I_{n}$ denote the $n\times n$-matrix
$\left(  \delta_{i,j}\right)  _{1\leq i\leq n,\ 1\leq j\leq n}$, where
$\delta_{i,j}$ is defined to be $%
\begin{cases}
1, & \text{if }i=j;\\
0, & \text{if }i\neq j
\end{cases}
$.\ \ \ \ \footnote{Here, $0$ and $1$ mean the zero and the unity of
$\mathbb{K}$ (which may and may not be the integers $0$ and $1$).} This matrix
$I_{n}$ looks as follows:%
\[
I_{n}=\left(
\begin{array}
[c]{cccc}%
1 & 0 & \cdots & 0\\
0 & 1 & \cdots & 0\\
\vdots & \vdots & \ddots & \vdots\\
0 & 0 & \cdots & 1
\end{array}
\right)  .
\]
It has the property that $I_{n}B=B$ for every $m\in\mathbb{N}$ and every
$n\times m$-matrix $B$; also, $AI_{n}=A$ for every $k\in\mathbb{N}$ and every
$k\times n$-matrix $A$. (Proving this is a good way to check that you
understand how matrices are multiplied.\footnote{See \cite[\S 2.12]{Gri-lina}
for a detailed proof of the equality $AI_{n}=A$. (Interpret the word
\textquotedblleft number\textquotedblright\ in \cite[\S 2.12]{Gri-lina} as
\textquotedblleft element of $\mathbb{K}$\textquotedblright.) The proof of
$I_{n}B=B$ is rather similar.}) The matrix $I_{n}$ is called the $n\times n$
\textit{identity matrix}. (Some call it $E_{n}$ or just $I$, when the value of
$n$ is clear from the context.)

Matrix multiplication is associative: If $n,m,k,\ell\in\mathbb{N}$, and if $A$
is an $n\times m$-matrix, $B$ is an $m\times k$-matrix, and $C$ is a
$k\times\ell$-matrix, then $A\left(  BC\right)  =\left(  AB\right)  C$. The
proof of this is straightforward using our definition of products of
matrices\footnote{Check that $A\left(  BC\right)  $ and $\left(  AB\right)  C$
both are equal to the matrix $\left(  \sum_{u=1}^{m}\sum_{v=1}^{k}%
a_{i,u}b_{u,v}c_{v,j}\right)  _{1\leq i\leq n,\ 1\leq j\leq\ell}$. For details
of this proof, see \cite[\S 2.9]{Gri-lina}. (Interpret the word
\textquotedblleft number\textquotedblright\ in \cite[\S 2.9]{Gri-lina} as
\textquotedblleft element of $\mathbb{K}$\textquotedblright.)}. This
associativity allows us to write products like $ABC$ without parentheses. By
induction, we can see that longer products such as $A_{1}A_{2}\cdots A_{k}$
for arbitrary $k\in\mathbb{N}$ can also be bracketed at will, because all
bracketings lead to the same result (e.g., for four matrices $A$, $B$, $C$ and
$D$, we have $A\left(  B\left(  CD\right)  \right)  =A\left(  \left(
BC\right)  D\right)  =\left(  AB\right)  \left(  CD\right)  =\left(  A\left(
BC\right)  \right)  D=\left(  \left(  AB\right)  C\right)  D$, provided that
the dimensions of the matrices are appropriate to make sense of the products).
We define an empty product of $n\times n$-matrices to be the $n\times n$
identity matrix $I_{n}$.

For every $n\times n$-matrix $A$ and every $k\in\mathbb{N}$, we can thus
define an $n\times n$-matrix $A^{k}$ by $A^{k}=\underbrace{AA\cdots
A}_{k\text{ factors}}$. In particular, $A^{0}=I_{n}$ (since we defined an
empty product of $n\times n$-matrices to be $I_{n}$).

Further properties of matrix multiplication are easy to state and to prove:

\begin{itemize}
\item For every $n\in\mathbb{N}$, $m\in\mathbb{N}$, $k\in\mathbb{N}$ and
$\lambda\in\mathbb{K}$, every $n\times m$-matrix $A$ and every $m\times
k$-matrix $B$, we have $\lambda\left(  AB\right)  =\left(  \lambda A\right)
B=A\left(  \lambda B\right)  $. (This allows us to write $\lambda AB$ for each
of the matrices $\lambda\left(  AB\right)  $, $\left(  \lambda A\right)  B$
and $A\left(  \lambda B\right)  $.)

\item For every $n\in\mathbb{N}$, $m\in\mathbb{N}$ and $k\in\mathbb{N}$, every
two $n\times m$-matrices $A$ and $B$, and every $m\times k$-matrix $C$, we
have $\left(  A+B\right)  C=AC+BC$.

\item For every $n\in\mathbb{N}$, $m\in\mathbb{N}$ and $k\in\mathbb{N}$, every
$n\times m$-matrix $A$, and every two $m\times k$-matrices $B$ and $C$, we
have $A\left(  B+C\right)  =AB+AC$.

\item For every $n\in\mathbb{N}$, $m\in\mathbb{N}$, $\lambda\in\mathbb{K}$ and
$\mu\in\mathbb{K}$, and every $n\times m$-matrix $A$, we have $\lambda\left(
\mu A\right)  =\left(  \lambda\mu\right)  A$. (This allows us to write
$\lambda\mu A$ for both $\lambda\left(  \mu A\right)  $ and $\left(
\lambda\mu\right)  A$.)
\end{itemize}

For given $n\in\mathbb{N}$ and $m\in\mathbb{N}$, we let $\mathbb{K}^{n\times
m}$ denote the set of all $n\times m$-matrices. (This is one of the two
standard notations for this set; the other is $\operatorname*{M}%
\nolimits_{n,m}\left(  \mathbb{K}\right)  $.)

A \textit{square matrix} is a matrix which has as many rows as it has columns;
in other words, a square matrix is an $n\times n$-matrix for some
$n\in\mathbb{N}$. If $A=\left(  a_{i,j}\right)  _{1\leq i\leq n,\ 1\leq j\leq
n}$ is a square matrix, then the $n$-tuple $\left(  a_{1,1},a_{2,2}%
,\ldots,a_{n,n}\right)  $ is called the \textit{diagonal} of $A$. (Some
authors abbreviate $\left(  a_{i,j}\right)  _{1\leq i\leq n,\ 1\leq j\leq n}$
by $\left(  a_{i,j}\right)  _{1\leq i,j\leq n}$; this notation has some mild
potential for confusion, though\footnote{The comma between \textquotedblleft%
$i$\textquotedblright\ and \textquotedblleft$j$\textquotedblright\ in
\textquotedblleft$1\leq i,j\leq n$\textquotedblright\ can be understood either
to separate $i$ from $j$, or to separate the inequality $1\leq i$ from the
inequality $j\leq n$. I remember seeing this ambiguity causing a real
misunderstanding.}.) The entries of the diagonal of $A$ are called the
\textit{diagonal entries} of $A$. (Some authors like to say \textquotedblleft
main diagonal\textquotedblright\ instead of \textquotedblleft
diagonal\textquotedblright.)

For a given $n\in\mathbb{N}$, the product of two $n\times n$-matrices is
always well-defined, and is an $n\times n$-matrix again. The set
$\mathbb{K}^{n\times n}$ satisfies all the axioms of a commutative ring except
for commutativity of multiplication. This makes it into what is commonly
called a \textit{noncommutative ring}\footnote{A \textit{noncommutative ring}
is defined in the same way as we defined a commutative ring, except for the
fact that commutativity of multiplication is removed from the list of axioms.
(The words \textquotedblleft noncommutative ring\textquotedblright\ do not
imply that commutativity of multiplication must be false for this ring; they
merely say that commutativity of multiplication is \textbf{not required} to
hold for it. For example, the noncommutative ring $\mathbb{K}^{n\times n}$ is
actually commutative when $n\leq1$ or when $\mathbb{K}$ is a trivial ring.)
\par
Instead of saying \textquotedblleft noncommutative ring\textquotedblright,
many algebraists just say \textquotedblleft ring\textquotedblright. We shall,
however, keep the word \textquotedblleft noncommutative\textquotedblright\ in
order to avoid confusion.}. We shall study noncommutative rings later (in
Section \ref{sect.noncommring}).

\subsection{Determinants}

Square matrices have determinants. Let us recall how determinants are defined:

\begin{definition}
\label{def.det}Let $n\in\mathbb{N}$. Let $A=\left(  a_{i,j}\right)  _{1\leq
i\leq n,\ 1\leq j\leq n}$ be an $n\times n$-matrix. The \textit{determinant}
$\det A$ of $A$ is defined as%
\begin{equation}
\sum_{\sigma\in S_{n}}\left(  -1\right)  ^{\sigma}a_{1,\sigma\left(  1\right)
}a_{2,\sigma\left(  2\right)  }\cdots a_{n,\sigma\left(  n\right)  }.
\label{eq.det}%
\end{equation}
In other words,%
\begin{align}
\det A  &  =\sum_{\sigma\in S_{n}}\left(  -1\right)  ^{\sigma}%
\underbrace{a_{1,\sigma\left(  1\right)  }a_{2,\sigma\left(  2\right)  }\cdots
a_{n,\sigma\left(  n\right)  }}_{=\prod_{i=1}^{n}a_{i,\sigma\left(  i\right)
}}\label{eq.det.eq.1}\\
&  =\sum_{\sigma\in S_{n}}\left(  -1\right)  ^{\sigma}\prod_{i=1}%
^{n}a_{i,\sigma\left(  i\right)  }. \label{eq.det.eq.2}%
\end{align}

\end{definition}

For example, the determinant of a $1\times1$-matrix $\left(
\begin{array}
[c]{c}%
a_{1,1}%
\end{array}
\right)  $ is%
\begin{align}
\det\left(
\begin{array}
[c]{c}%
a_{1,1}%
\end{array}
\right)   &  =\sum_{\sigma\in S_{1}}\left(  -1\right)  ^{\sigma}%
a_{1,\sigma\left(  1\right)  }=\underbrace{\left(  -1\right)
^{\operatorname*{id}}}_{=1}\underbrace{a_{1,\operatorname*{id}\left(
1\right)  }}_{=a_{1,1}}\nonumber\\
&  \ \ \ \ \ \ \ \ \ \ \ \ \ \ \ \ \ \ \ \ \left(  \text{since the only
permutation }\sigma\in S_{1}\text{ is }\operatorname*{id}\right) \nonumber\\
&  =a_{1,1}. \label{eq.det.small.1x1}%
\end{align}
The determinant of a $2\times2$-matrix $\left(
\begin{array}
[c]{cc}%
a_{1,1} & a_{1,2}\\
a_{2,1} & a_{2,2}%
\end{array}
\right)  $ is%
\begin{align*}
\det\left(
\begin{array}
[c]{cc}%
a_{1,1} & a_{1,2}\\
a_{2,1} & a_{2,2}%
\end{array}
\right)   &  =\sum_{\sigma\in S_{2}}\left(  -1\right)  ^{\sigma}%
a_{1,\sigma\left(  1\right)  }a_{2,\sigma\left(  2\right)  }\\
&  =\underbrace{\left(  -1\right)  ^{\operatorname*{id}}}_{=1}%
\underbrace{a_{1,\operatorname*{id}\left(  1\right)  }}_{=a_{1,1}%
}\underbrace{a_{2,\operatorname*{id}\left(  2\right)  }}_{=a_{2,2}%
}+\underbrace{\left(  -1\right)  ^{s_{1}}}_{=-1}\underbrace{a_{1,s_{1}\left(
1\right)  }}_{=a_{1,2}}\underbrace{a_{2,s_{1}\left(  2\right)  }}_{=a_{2,1}}\\
&  \ \ \ \ \ \ \ \ \ \ \ \ \ \ \ \ \ \ \ \ \left(  \text{since the only
permutations }\sigma\in S_{2}\text{ are }\operatorname*{id}\text{ and }%
s_{1}\right) \\
&  =a_{1,1}a_{2,2}-a_{1,2}a_{2,1}.
\end{align*}
Similarly, for a $3\times3$-matrix, the formula is%
\begin{align}
\det\left(
\begin{array}
[c]{ccc}%
a_{1,1} & a_{1,2} & a_{1,3}\\
a_{2,1} & a_{2,2} & a_{2,3}\\
a_{3,1} & a_{3,2} & a_{3,3}%
\end{array}
\right)   &  =a_{1,1}a_{2,2}a_{3,3}+a_{1,2}a_{2,3}a_{3,1}+a_{1,3}%
a_{2,1}a_{3,2}\nonumber\\
&  \ \ \ \ \ \ \ \ \ \ -a_{1,1}a_{2,3}a_{3,2}-a_{1,2}a_{2,1}a_{3,3}%
-a_{1,3}a_{2,2}a_{3,1}. \label{eq.det.small.3x3}%
\end{align}
Also, the determinant of the $0\times0$-matrix is $1$\ \ \ \ \footnote{In more
details:
\par
There is only one $0\times0$-matrix; it has no rows and no columns and no
entries. According to (\ref{eq.det.eq.2}), its determinant is
\begin{align*}
\sum_{\sigma\in S_{0}}\left(  -1\right)  ^{\sigma}\underbrace{\prod_{i=1}%
^{0}a_{i,\sigma\left(  i\right)  }}_{=\left(  \text{empty product}\right)
=1}  &  =\sum_{\sigma\in S_{0}}\left(  -1\right)  ^{\sigma}=\left(  -1\right)
^{\operatorname*{id}}\ \ \ \ \ \ \ \ \ \ \left(  \text{since the only }%
\sigma\in S_{0}\text{ is }\operatorname*{id}\right) \\
&  =1.
\end{align*}
}. (This might sound like hairsplitting, but being able to work with
$0\times0$-matrices simplifies some proofs by induction, because it allows one
to take $n=0$ as an induction base.)

The equality (\ref{eq.det.eq.2}) (or, equivalently, (\ref{eq.det.eq.1})) is
known as the \textit{Leibniz formula}. Out of several known ways to define the
determinant, it is probably the most direct. In practice, however, computing a
determinant using (\ref{eq.det.eq.2}) quickly becomes impractical when $n$ is
high (since the sum has $n!$ terms). In most situations that occur both
\href{http://arxiv.org/abs/math/9902004v3}{in mathematics} and
\href{https://en.wikipedia.org/wiki/Determinant#Calculation}{in applications},
determinants can be computed in various simpler ways.

Some authors write $\left\vert A\right\vert $ instead of $\det A$ for the
determinant of a square matrix $A$. I do not like this notation, as it clashes
(in the case of $1\times1$-matrices) with the notation $\left\vert
a\right\vert $ for the absolute value of a real number $a$.

Here is a first example of a determinant which ends up very simple:

\begin{example}
\label{exam.xiyj}Let $n\in\mathbb{N}$. Let $x_{1},x_{2},\ldots,x_{n}$ be $n$
elements of $\mathbb{K}$, and let $y_{1},y_{2},\ldots,y_{n}$ be $n$ further
elements of $\mathbb{K}$. Let $A$ be the $n\times n$-matrix $\left(
x_{i}y_{j}\right)  _{1\leq i\leq n,\ 1\leq j\leq n}$. What is $\det A$ ?

For $n=0$, we have $\det A=1$ (since the $0\times0$-matrix has determinant $1$).

For $n=1$, we have $A=\left(
\begin{array}
[c]{c}%
x_{1}y_{1}%
\end{array}
\right)  $ and thus $\det A=x_{1}y_{1}$.

For $n=2$, we have $A=\left(
\begin{array}
[c]{cc}%
x_{1}y_{1} & x_{1}y_{2}\\
x_{2}y_{1} & x_{2}y_{2}%
\end{array}
\right)  $ and thus $\det A=\left(  x_{1}y_{1}\right)  \left(  x_{2}%
y_{2}\right)  -\left(  x_{1}y_{2}\right)  \left(  x_{2}y_{1}\right)  =0$.

What do you expect for greater values of $n$ ? The pattern might not be clear
at this point yet, but if you compute further examples, you will realize that
$\det A=0$ also holds for $n=3$, for $n=4$, for $n=5$... This suggests that
$\det A=0$ for every $n\geq2$. How to prove this?

Let $n\geq2$. Then, (\ref{eq.det.eq.1}) (applied to $a_{i,j}=x_{i}y_{j}$)
yields%
\begin{align}
\det A  &  =\sum_{\sigma\in S_{n}}\left(  -1\right)  ^{\sigma}%
\underbrace{\left(  x_{1}y_{\sigma\left(  1\right)  }\right)  \left(
x_{2}y_{\sigma\left(  2\right)  }\right)  \cdots\left(  x_{n}y_{\sigma\left(
n\right)  }\right)  }_{=\left(  x_{1}x_{2}\cdots x_{n}\right)  \left(
y_{\sigma\left(  1\right)  }y_{\sigma\left(  2\right)  }\cdots y_{\sigma
\left(  n\right)  }\right)  }\nonumber\\
&  =\sum_{\sigma\in S_{n}}\left(  -1\right)  ^{\sigma}\left(  x_{1}x_{2}\cdots
x_{n}\right)  \underbrace{\left(  y_{\sigma\left(  1\right)  }y_{\sigma\left(
2\right)  }\cdots y_{\sigma\left(  n\right)  }\right)  }_{\substack{=y_{1}%
y_{2}\cdots y_{n}\\\text{(since }\sigma\text{ is a permutation)}}}\nonumber\\
&  =\sum_{\sigma\in S_{n}}\left(  -1\right)  ^{\sigma}\left(  x_{1}x_{2}\cdots
x_{n}\right)  \left(  y_{1}y_{2}\cdots y_{n}\right) \nonumber\\
&  =\left(  \sum_{\sigma\in S_{n}}\left(  -1\right)  ^{\sigma}\right)  \left(
x_{1}x_{2}\cdots x_{n}\right)  \left(  y_{1}y_{2}\cdots y_{n}\right)  .
\label{eq.exam.xiyj.detA1}%
\end{align}
Now, every $\sigma\in S_{n}$ is either even or odd (but not both), and thus we
have%
\begin{align*}
&  \sum_{\sigma\in S_{n}}\left(  -1\right)  ^{\sigma}\\
&  =\sum_{\substack{\sigma\in S_{n};\\\sigma\text{ is even}}%
}\underbrace{\left(  -1\right)  ^{\sigma}}_{\substack{=1\\\text{(since }%
\sigma\text{ is even)}}}+\sum_{\substack{\sigma\in S_{n};\\\sigma\text{ is
odd}}}\underbrace{\left(  -1\right)  ^{\sigma}}_{\substack{=-1\\\text{(since
}\sigma\text{ is odd)}}}\\
&  =\underbrace{\sum_{\substack{\sigma\in S_{n};\\\sigma\text{ is even}}%
}1}_{=\left(  \text{the number of even permutations }\sigma\in S_{n}\right)
\cdot1}+\underbrace{\sum_{\substack{\sigma\in S_{n};\\\sigma\text{ is odd}%
}}\left(  -1\right)  }_{=\left(  \text{the number of odd permutations }%
\sigma\in S_{n}\right)  \cdot\left(  -1\right)  }\\
&  =\underbrace{\left(  \text{the number of even permutations }\sigma\in
S_{n}\right)  }_{\substack{=n!/2\\\text{(by Exercise \ref{exe.ps2.2.7})}%
}}\cdot1\\
&  \ \ \ \ \ \ \ \ \ \ +\underbrace{\left(  \text{the number of odd
permutations }\sigma\in S_{n}\right)  }_{\substack{=n!/2\\\text{(by Exercise
\ref{exe.ps2.2.7})}}}\cdot\left(  -1\right) \\
&  =\left(  n!/2\right)  \cdot1+\left(  n!/2\right)  \cdot\left(  -1\right)
=0.
\end{align*}
Hence, (\ref{eq.exam.xiyj.detA1}) becomes $\det A=\underbrace{\left(
\sum_{\sigma\in S_{n}}\left(  -1\right)  ^{\sigma}\right)  }_{=0}\left(
x_{1}x_{2}\cdots x_{n}\right)  \left(  y_{1}y_{2}\cdots y_{n}\right)  =0$, as
we wanted to prove.

We will eventually learn a simpler way to prove this.
\end{example}

\begin{example}
\label{exam.xi+yj}Here is an example similar to Example \ref{exam.xiyj}, but subtler.

Let $n\in\mathbb{N}$. Let $x_{1},x_{2},\ldots,x_{n}$ be $n$ elements of
$\mathbb{K}$, and let $y_{1},y_{2},\ldots,y_{n}$ be $n$ further elements of
$\mathbb{K}$. Let $A$ be the $n\times n$-matrix $\left(  x_{i}+y_{j}\right)
_{1\leq i\leq n,\ 1\leq j\leq n}$. What is $\det A$ ?

For $n=0$, we have $\det A=1$ again.

For $n=1$, we have $A=\left(
\begin{array}
[c]{c}%
x_{1}+y_{1}%
\end{array}
\right)  $ and thus $\det A=x_{1}+y_{1}$.

For $n=2$, we have $A=\left(
\begin{array}
[c]{cc}%
x_{1}+y_{1} & x_{1}+y_{2}\\
x_{2}+y_{1} & x_{2}+y_{2}%
\end{array}
\right)  $ and thus $\det A=\left(  x_{1}+y_{1}\right)  \left(  x_{2}%
+y_{2}\right)  -\left(  x_{1}+y_{2}\right)  \left(  x_{2}+y_{1}\right)
=-\left(  y_{1}-y_{2}\right)  \left(  x_{1}-x_{2}\right)  $.

However, it turns out that for every $n\geq3$, we again have $\det A=0$. This
is harder to prove than the similar claim in Example \ref{exam.xiyj}. We will
eventually see how to do it easily, but as for now let me outline a direct
proof. (I am being rather telegraphic here; do not worry if you do not
understand the following argument, as there will be easier and more detailed
proofs below.)

From (\ref{eq.det.eq.1}), we obtain%
\begin{equation}
\det A=\sum_{\sigma\in S_{n}}\left(  -1\right)  ^{\sigma}\left(
x_{1}+y_{\sigma\left(  1\right)  }\right)  \left(  x_{2}+y_{\sigma\left(
2\right)  }\right)  \cdots\left(  x_{n}+y_{\sigma\left(  n\right)  }\right)  .
\label{eq.exam.xi+yj.detA1}%
\end{equation}
If we expand the product $\left(  x_{1}+y_{\sigma\left(  1\right)  }\right)
\left(  x_{2}+y_{\sigma\left(  2\right)  }\right)  \cdots\left(
x_{n}+y_{\sigma\left(  n\right)  }\right)  $, we obtain a sum of $2^{n}$
terms:%
\begin{align}
&  \left(  x_{1}+y_{\sigma\left(  1\right)  }\right)  \left(  x_{2}%
+y_{\sigma\left(  2\right)  }\right)  \cdots\left(  x_{n}+y_{\sigma\left(
n\right)  }\right) \nonumber\\
&  =\sum_{I\subseteq\left[  n\right]  }\left(  \prod_{i\in I}x_{i}\right)
\left(  \prod_{i\in\left[  n\right]  \setminus I}y_{\sigma\left(  i\right)
}\right)  \label{eq.exam.xi+yj.detA2}%
\end{align}
(where $\left[  n\right]  $ means the set $\left\{  1,2,\ldots,n\right\}  $).
(To obtain a fully rigorous proof of (\ref{eq.exam.xi+yj.detA2}), apply
Exercise \ref{exe.prod(ai+bi)} \textbf{(a)} to $a_{i}=x_{i}$ and
$b_{i}=y_{\sigma\left(  i\right)  }$.) Thus, (\ref{eq.exam.xi+yj.detA1})
becomes%
\begin{align*}
\det A  &  =\sum_{\sigma\in S_{n}}\left(  -1\right)  ^{\sigma}%
\underbrace{\left(  x_{1}+y_{\sigma\left(  1\right)  }\right)  \left(
x_{2}+y_{\sigma\left(  2\right)  }\right)  \cdots\left(  x_{n}+y_{\sigma
\left(  n\right)  }\right)  }_{\substack{=\sum_{I\subseteq\left[  n\right]
}\left(  \prod_{i\in I}x_{i}\right)  \left(  \prod_{i\in\left[  n\right]
\setminus I}y_{\sigma\left(  i\right)  }\right)  \\\text{(by
(\ref{eq.exam.xi+yj.detA2}))}}}\\
&  =\sum_{\sigma\in S_{n}}\left(  -1\right)  ^{\sigma}\sum_{I\subseteq\left[
n\right]  }\left(  \prod_{i\in I}x_{i}\right)  \left(  \prod_{i\in\left[
n\right]  \setminus I}y_{\sigma\left(  i\right)  }\right) \\
&  =\sum_{I\subseteq\left[  n\right]  }\sum_{\sigma\in S_{n}}\left(
-1\right)  ^{\sigma}\left(  \prod_{i\in I}x_{i}\right)  \left(  \prod
_{i\in\left[  n\right]  \setminus I}y_{\sigma\left(  i\right)  }\right)  .
\end{align*}
We want to prove that this is $0$. In order to do so, it clearly suffices to
show that every $I\subseteq\left[  n\right]  $ satisfies%
\begin{equation}
\sum_{\sigma\in S_{n}}\left(  -1\right)  ^{\sigma}\left(  \prod_{i\in I}%
x_{i}\right)  \left(  \prod_{i\in\left[  n\right]  \setminus I}y_{\sigma
\left(  i\right)  }\right)  =0. \label{eq.exam.xi+yj.detA4}%
\end{equation}
So let us fix $I\subseteq\left[  n\right]  $, and try to prove
(\ref{eq.exam.xi+yj.detA4}). We must be in one of the following two cases:

\begin{description}
\item[Case 1:] The set $\left[  n\right]  \setminus I$ has at least two
elements. In this case, let us pick two distinct elements $a$ and $b$ of this
set, and split the set $S_{n}$ into disjoint two-element subsets by pairing up
every even permutation $\sigma\in S_{n}$ with the odd permutation $\sigma\circ
t_{a,b}$ (where $t_{a,b}$ is as defined in Definition \ref{def.transpos}). The
addends on the left hand side of (\ref{eq.exam.xi+yj.detA4}) corresponding to
two permutations paired up cancel out each other (because the products
$\prod_{i\in\left[  n\right]  \setminus I}y_{\sigma\left(  i\right)  }$ and
$\prod_{i\in\left[  n\right]  \setminus I}y_{\left(  \sigma\circ
t_{a,b}\right)  \left(  i\right)  }$ differ only in the order of their
factors), and thus the whole left hand side of (\ref{eq.exam.xi+yj.detA4}) is
$0$. Thus, (\ref{eq.exam.xi+yj.detA4}) is proven in Case 1.

\item[Case 2:] The set $\left[  n\right]  \setminus I$ has at most one
element. In this case, the set $I$ has at least two elements (it is here that
we use $n\geq3$). Pick two distinct elements $c$ and $d$ of $I$, and split the
set $S_{n}$ into disjoint two-element subsets by pairing up every even
permutation $\sigma\in S_{n}$ with the odd permutation $\sigma\circ t_{c,d}$.
Again, the addends on the left hand side of (\ref{eq.exam.xi+yj.detA4})
corresponding to two permutations paired up cancel out each other (because the
products $\prod_{i\in\left[  n\right]  \setminus I}y_{\sigma\left(  i\right)
}$ and $\prod_{i\in\left[  n\right]  \setminus I}y_{\left(  \sigma\circ
t_{c,d}\right)  \left(  i\right)  }$ are identical), and thus the whole left
hand side of (\ref{eq.exam.xi+yj.detA4}) is $0$. This proves
(\ref{eq.exam.xi+yj.detA4}) in Case 2.
\end{description}

We thus have proven (\ref{eq.exam.xi+yj.detA4}) in both cases. So $\det A=0$
is proven. This was a tricky argument, and shows the limits of the usefulness
of (\ref{eq.det.eq.1}).
\end{example}

We shall now discuss basic properties of the determinant.

\begin{exercise}
\label{exe.ps4.3}Let $A=\left(  a_{i,j}\right)  _{1\leq i\leq n,\ 1\leq j\leq
n}$ be an $n\times n$-matrix. Assume that $a_{i,j}=0$ for every $\left(
i,j\right)  \in\left\{  1,2,\ldots,n\right\}  ^{2}$ satisfying $i<j$. Show
that%
\[
\det A=a_{1,1}a_{2,2}\cdots a_{n,n}.
\]

\end{exercise}

\begin{definition}
An $n\times n$-matrix $A$ satisfying the assumption of Exercise
\ref{exe.ps4.3} is said to be \textit{lower-triangular} (because its entries
above the diagonal are $0$, and thus its nonzero entries are concentrated in
the triangular region southwest of the diagonal). Exercise \ref{exe.ps4.3}
thus says that the determinant of a lower-triangular matrix is the product of
its diagonal entries. For instance, $\det\left(
\begin{array}
[c]{ccc}%
a & 0 & 0\\
b & c & 0\\
d & e & f
\end{array}
\right)  =acf$.
\end{definition}

\begin{example}
\label{exa.det.In}Let $n\in\mathbb{N}$. The $n\times n$ identity matrix
$I_{n}$ is lower-triangular, and its diagonal entries are $1,1,\ldots,1$.
Hence, Exercise \ref{exe.ps4.3} shows that its determinant is $\det\left(
I_{n}\right)  =1\cdot1\cdot\cdots\cdot1=1$.
\end{example}

\begin{definition}
\label{def.transpose}The \textit{transpose} of a matrix $A=\left(
a_{i,j}\right)  _{1\leq i\leq n,\ 1\leq j\leq m}$ is defined to be the matrix
$\left(  a_{j,i}\right)  _{1\leq i\leq m,\ 1\leq j\leq n}$. It is denoted by
$A^{T}$. For instance, $\left(
\begin{array}
[c]{ccc}%
1 & 2 & -1\\
4 & 0 & 1
\end{array}
\right)  ^{T}=\left(
\begin{array}
[c]{cc}%
1 & 4\\
2 & 0\\
-1 & 1
\end{array}
\right)  $.
\end{definition}

\begin{remark}
Various other notations for the transpose of a matrix $A$ exist in the
literature. Some of them are $A^{t}$ (with a lower case $t$) and $^{T}A$ and
$^{t}A$.
\end{remark}

\begin{exercise}
\label{exe.ps4.4}Let $n\in\mathbb{N}$. Let $A$ be an $n\times n$-matrix. Show
that $\det\left(  A^{T}\right)  =\det A$.
\end{exercise}

The transpose of a lower-triangular $n\times n$-matrix is an upper-triangular
$n\times n$-matrix (i.e., an $n\times n$-matrix whose entries below the
diagonal are $0$). Thus, combining Exercise \ref{exe.ps4.3} with Exercise
\ref{exe.ps4.4}, we see that the determinant of an upper-triangular matrix is
the product of its diagonal entries.

The following exercise presents five fundamental (and simple) properties of transposes:

\begin{exercise}
\label{exe.transpose.basics}Prove the following:

\textbf{(a)} If $u$, $v$ and $w$ are three nonnegative integers, if $P$ is a
$u\times v$-matrix, and if $Q$ is a $v\times w$-matrix, then%
\begin{equation}
\left(  PQ\right)  ^{T}=Q^{T}P^{T}. \label{pf.thm.adj.inverse.tranposes1}%
\end{equation}

\textbf{(b)} Every $u\in\mathbb{N}$ satisfies%
\begin{equation}
\left(  I_{u}\right)  ^{T}=I_{u}. \label{pf.thm.adj.inverse.tranposes2}%
\end{equation}

\textbf{(c)} If $u$ and $v$ are two nonnegative integers, if $P$ is a $u\times
v$-matrix, and if $\lambda\in\mathbb{K}$, then%
\begin{equation}
\left(  \lambda P\right)  ^{T}=\lambda P^{T}.
\label{pf.thm.adj.inverse.tranposes3}%
\end{equation}

\textbf{(d)} If $u$ and $v$ are two nonnegative integers, and if $P$ and $Q$
are two $u\times v$-matrices, then%
\[
\left(  P+Q\right)  ^{T}=P^{T}+Q^{T}.
\]

\textbf{(e)} If $u$ and $v$ are two nonnegative integers, and if $P$ is a
$u\times v$-matrix, then%
\begin{equation}
\left(  P^{T}\right)  ^{T}=P. \label{pf.thm.adj.inverse.tranposes4}%
\end{equation}

\end{exercise}

Here is yet another simple property of determinants that follows directly from
their definition:

\begin{proposition}
\label{prop.det.scale}Let $n\in\mathbb{N}$ and $\lambda\in\mathbb{K}$. Let $A$
be an $n\times n$-matrix. Then, $\det\left(  \lambda A\right)  =\lambda
^{n}\det A$.
\end{proposition}

\begin{proof}
[Proof of Proposition \ref{prop.det.scale}.]Write $A$ in the form $A=\left(
a_{i,j}\right)  _{1\leq i\leq n,\ 1\leq j\leq n}$. Thus, $\lambda A=\left(
\lambda a_{i,j}\right)  _{1\leq i\leq n,\ 1\leq j\leq n}$ (by the definition
of $\lambda A$). Hence, (\ref{eq.det.eq.2}) (applied to $\lambda A$ and
$\lambda a_{i,j}$ instead of $A$ and $a_{i,j}$) yields%
\begin{align*}
\det\left(  \lambda A\right)   &  =\sum_{\sigma\in S_{n}}\left(  -1\right)
^{\sigma}\underbrace{\prod_{i=1}^{n}\left(  \lambda a_{i,\sigma\left(
i\right)  }\right)  }_{=\lambda^{n}\prod_{i=1}^{n}a_{i,\sigma\left(  i\right)
}}=\sum_{\sigma\in S_{n}}\left(  -1\right)  ^{\sigma}\lambda^{n}\prod
_{i=1}^{n}a_{i,\sigma\left(  i\right)  }\\
&  =\lambda^{n}\underbrace{\sum_{\sigma\in S_{n}}\left(  -1\right)  ^{\sigma
}\prod_{i=1}^{n}a_{i,\sigma\left(  i\right)  }}_{\substack{=\det A\\\text{(by
(\ref{eq.det.eq.2}))}}}=\lambda^{n}\det A.
\end{align*}
Proposition \ref{prop.det.scale} is thus proven.
\end{proof}

\Needspace{8cm}

\begin{exercise}
\label{exe.ps4.5}Let $a,b,c,d,e,f,g,h,i,j,k,\ell,m,n,o,p$ be elements of
$\mathbb{K}$.

\textbf{(a)} Find a simple formula for the determinant%
\[
\det\left(
\begin{array}
[c]{cccc}%
a & b & c & d\\
\ell & 0 & 0 & e\\
k & 0 & 0 & f\\
j & i & h & g
\end{array}
\right)  .
\]

\textbf{(b)} Find a simple formula for the determinant%
\[
\det\left(
\begin{array}
[c]{ccccc}%
a & b & c & d & e\\
f & 0 & 0 & 0 & g\\
h & 0 & 0 & 0 & i\\
j & 0 & 0 & 0 & k\\
\ell & m & n & o & p
\end{array}
\right)  .
\]
(Do not mistake the \textquotedblleft$o$\textquotedblright\ for a
\textquotedblleft$0$\textquotedblright.)

[\textbf{Hint:} Part \textbf{(b)} is simpler than part \textbf{(a)}.]
\end{exercise}

In the next exercises, we shall talk about rows and columns; let us first make
some pedantic remarks about these notions.

If $n\in\mathbb{N}$, then an $n\times1$-matrix is said to be a \textit{column
vector} with $n$ entries\footnote{It is also called a \textit{column vector}
of size $n$.}, whereas a $1\times n$-matrix is said to be a \textit{row
vector} with $n$ entries. Column vectors and row vectors store exactly the
same kind of data (namely, $n$ elements of $\mathbb{K}$), so you might wonder
why I make a difference between them (and also why I distinguish them from
$n$-tuples of elements of $\mathbb{K}$, which also contain precisely the same
kind of data). The reason for this is that column vectors and row vectors
behave differently under matrix multiplication: For example,%
\[
\left(
\begin{array}
[c]{c}%
a\\
b
\end{array}
\right)  \left(
\begin{array}
[c]{cc}%
c & d
\end{array}
\right)  =\left(
\begin{array}
[c]{cc}%
ac & ad\\
bc & bd
\end{array}
\right)
\]
is not the same as%
\[
\left(
\begin{array}
[c]{cc}%
a & b
\end{array}
\right)  \left(
\begin{array}
[c]{c}%
c\\
d
\end{array}
\right)  =\left(
\begin{array}
[c]{c}%
ac+bd
\end{array}
\right)  .
\]
If we would identify column vectors with row vectors, then this would cause contradictions.

The reason to distinguish between row vectors and $n$-tuples is subtler: We
have defined row vectors only for a commutative ring $\mathbb{K}$, whereas
$n$-tuples can be made out of elements of any set. As a consequence, the sum
of two row vectors is well-defined (since row vectors are matrices and thus
can be added entry by entry), whereas the sum of two $n$-tuples is not.
Similarly, we can take the product $\lambda v$ of an element $\lambda
\in\mathbb{K}$ with a row vector $v$ (by multiplying every entry of $v$ by
$\lambda$), but such a thing does not make sense for general $n$-tuples. These
differences between row vectors and $n$-tuples, however, cause no clashes of
notation if we use the same notations for both types of object. Thus, we are
often going to identify a row vector $\left(
\begin{array}
[c]{cccc}%
a_{1} & a_{2} & \cdots & a_{n}%
\end{array}
\right)  $ with the $n$-tuple $\left(  a_{1},a_{2},\ldots,a_{n}\right)
\in\mathbb{K}^{n}$. Thus, $\mathbb{K}^{n}$ becomes the set of all row vectors
with $n$ entries.\footnote{Some algebraists, instead, identify column vectors
with $n$-tuples, so that $\mathbb{K}^{n}$ is then the set of all column
vectors with $n$ entries. This is a valid convention as well, but one must be
careful not to use it simultaneously with our convention (i.e., with the
convention that row vectors are identified with $n$-tuples); this is why we
will not use it.}

The column vectors with $n$ entries are in 1-to-1 correspondence with the row
vectors with $n$ entries, and this correspondence is given by taking the
transpose: The column vector $v$ corresponds to the row vector $v^{T}$, and
conversely, the row vector $w$ corresponds to the column vector $w^{T}$. In
particular, every column vector $\left(
\begin{array}
[c]{c}%
a_{1}\\
a_{2}\\
\vdots\\
a_{n}%
\end{array}
\right)  $ can be rewritten in the form $\left(
\begin{array}
[c]{cccc}%
a_{1} & a_{2} & \cdots & a_{n}%
\end{array}
\right)  ^{T}=\left(  a_{1},a_{2},\ldots,a_{n}\right)  ^{T}$. We shall often
write it in the latter form, just because it takes up less space on paper.

The rows of a matrix are row vectors; the columns of a matrix are column
vectors. Thus, terms like \textquotedblleft the sum of two rows of a matrix
$A$\textquotedblright\ or \textquotedblleft$-3$ times a column of a matrix
$A$\textquotedblright\ make sense: Rows and columns are vectors, and thus can
be added (when they have the same number of entries) and multiplied by
elements of $\mathbb{K}$.

Let $n\in\mathbb{N}$ and $j\in\left\{  1,2,\ldots,n\right\}  $. If $v$ is a
column vector with $n$ entries (that is, an $n\times1$-matrix), then the
$j$\textit{-th entry of }$v$ means the $\left(  j,1\right)  $-th entry of $v$.
If $v$ is a row vector with $n$ entries (that is, a $1\times n$-matrix), then
the $j$\textit{-th entry of }$v$ means the $\left(  1,j\right)  $-th entry of
$v$. For example, the $2$-nd entry of the row vector $\left(
\begin{array}
[c]{ccc}%
a & b & c
\end{array}
\right)  $ is $b$.

\begin{exercise}
\label{exe.ps4.6}Let $n\in\mathbb{N}$. Let $A$ be an $n\times n$-matrix. Prove
the following:

\textbf{(a)} If $B$ is an $n\times n$-matrix obtained from $A$ by swapping two
rows, then $\det B=-\det A$. (\textquotedblleft Swapping two
rows\textquotedblright\ means \textquotedblleft swapping two distinct
rows\textquotedblright, of course.)

\textbf{(b)} If $B$ is an $n\times n$-matrix obtained from $A$ by swapping two
columns, then $\det B=-\det A$.

\textbf{(c)} If a row of $A$ consists of zeroes, then $\det A=0$.

\textbf{(d)} If a column of $A$ consists of zeroes, then $\det A=0$.

\textbf{(e)} If $A$ has two equal rows, then $\det A=0$.

\textbf{(f)} If $A$ has two equal columns, then $\det A=0$.

\textbf{(g)} Let $\lambda\in\mathbb{K}$ and $k\in\left\{  1,2,\ldots
,n\right\}  $. If $B$ is the $n\times n$-matrix obtained from $A$ by
multiplying the $k$-th row by $\lambda$ (that is, multiplying every entry of
the $k$-th row by $\lambda$), then $\det B=\lambda\det A$.

\textbf{(h)} Let $\lambda\in\mathbb{K}$ and $k\in\left\{  1,2,\ldots
,n\right\}  $. If $B$ is the $n\times n$-matrix obtained from $A$ by
multiplying the $k$-th column by $\lambda$, then $\det B=\lambda\det A$.

\textbf{(i)} Let $k\in\left\{  1,2,\ldots,n\right\}  $. Let $A^{\prime}$ be an
$n\times n$-matrix whose rows equal the corresponding rows of $A$ except
(perhaps) the $k$-th row. Let $B$ be the $n\times n$-matrix obtained from $A$
by adding the $k$-th row of $A^{\prime}$ to the $k$-th row of $A$ (that is, by
adding every entry of the $k$-th row of $A^{\prime}$ to the corresponding
entry of the $k$-th row of $A$). Then, $\det B=\det A+\det A^{\prime}$.

\textbf{(j)} Let $k\in\left\{  1,2,\ldots,n\right\}  $. Let $A^{\prime}$ be an
$n\times n$-matrix whose columns equal the corresponding columns of $A$ except
(perhaps) the $k$-th column. Let $B$ be the $n\times n$-matrix obtained from
$A$ by adding the $k$-th column of $A^{\prime}$ to the $k$-th column of $A$.
Then, $\det B=\det A+\det A^{\prime}$.
\end{exercise}

\begin{example}
Let us show examples for several parts of Exercise \ref{exe.ps4.6}
(especially, for Exercise \ref{exe.ps4.6} \textbf{(i)}, which has a somewhat
daunting statement).

\textbf{(a)} Exercise \ref{exe.ps4.6} \textbf{(a)} yields (among other things)
that
\[
\det\left(
\begin{array}
[c]{ccc}%
a & b & c\\
d & e & f\\
g & h & i
\end{array}
\right)  =-\det\left(
\begin{array}
[c]{ccc}%
g & h & i\\
d & e & f\\
a & b & c
\end{array}
\right)
\]
for any $a,b,c,d,e,f,g,h,i\in\mathbb{K}$.

\textbf{(c)} Exercise \ref{exe.ps4.6} \textbf{(c)} yields (among other things)
that
\[
\det\left(
\begin{array}
[c]{ccc}%
a & b & c\\
0 & 0 & 0\\
d & e & f
\end{array}
\right)  =0
\]
for any $a,b,c,d,e,f\in\mathbb{K}$.

\textbf{(e)} Exercise \ref{exe.ps4.6} \textbf{(e)} yields (among other things)
that
\[
\det\left(
\begin{array}
[c]{ccc}%
a & b & c\\
d & e & f\\
d & e & f
\end{array}
\right)  =0
\]
for any $a,b,c,d,e,f\in\mathbb{K}$.

\textbf{(g)} Exercise \ref{exe.ps4.6} \textbf{(g)} (applied to $n=3$ and
$k=2$) yields that
\[
\det\left(
\begin{array}
[c]{ccc}%
a & b & c\\
\lambda d & \lambda e & \lambda f\\
g & h & i
\end{array}
\right)  =\lambda\det\left(
\begin{array}
[c]{ccc}%
a & b & c\\
d & e & f\\
g & h & i
\end{array}
\right)
\]
for any $a,b,c,d,e,f\in\mathbb{K}$ and $\lambda\in\mathbb{K}$.

\textbf{(i)} Set $n=3$ and $k=2$. Set $A=\left(
\begin{array}
[c]{ccc}%
a & b & c\\
d & e & f\\
g & h & i
\end{array}
\right)  $. Then, a matrix $A^{\prime}$ satisfying the conditions of Exercise
\ref{exe.ps4.6} \textbf{(i)} has the form $A^{\prime}=\left(
\begin{array}
[c]{ccc}%
a & b & c\\
d^{\prime} & e^{\prime} & f^{\prime}\\
g & h & i
\end{array}
\right)  $. For such a matrix $A^{\prime}$, we obtain $B=\left(
\begin{array}
[c]{ccc}%
a & b & c\\
d+d^{\prime} & e+e^{\prime} & f+f^{\prime}\\
g & h & i
\end{array}
\right)  $. Exercise \ref{exe.ps4.6} \textbf{(i)} then claims that $\det
B=\det A+\det A^{\prime}$.
\end{example}

Parts \textbf{(a)}, \textbf{(c)}, \textbf{(e)}, \textbf{(g)} and \textbf{(i)}
of Exercise \ref{exe.ps4.6} are often united under the slogan
\textquotedblleft the determinant of a matrix is multilinear and alternating
in its rows\textquotedblright\footnote{Specifically, parts \textbf{(c)},
\textbf{(g)} and \textbf{(i)} say that it is \textquotedblleft
multilinear\textquotedblright, while parts \textbf{(a)} and \textbf{(e)} are
responsible for the \textquotedblleft alternating\textquotedblright.}.
Similarly, parts \textbf{(b)}, \textbf{(d)}, \textbf{(f)}, \textbf{(h)} and
\textbf{(j)} are combined under the slogan \textquotedblleft the determinant
of a matrix is multilinear and alternating in its columns\textquotedblright.
Many texts on linear algebra (for example, \cite{HoffmanKunze}) use these
properties as the \textbf{definition} of the determinant\footnote{More
precisely, they define a \textit{determinant function} to be a function
$F:\mathbb{K}^{n\times n}\rightarrow\mathbb{K}$ which is multilinear and
alternating in the rows of a matrix (i.e., which satisfies parts \textbf{(a)},
\textbf{(c)}, \textbf{(e)}, \textbf{(g)} and \textbf{(i)} of Exercise
\ref{exe.ps4.6} if every appearance of \textquotedblleft$\det$%
\textquotedblright\ is replaced by \textquotedblleft$F$\textquotedblright\ in
this Exercise) and which satisfies $F\left(  I_{n}\right)  =1$. Then, they
show that there is (for each $n\in\mathbb{N}$) exactly one determinant
function $F:\mathbb{K}^{n\times n}\rightarrow\mathbb{K}$. They then denote
this function by $\det$. This is a rather slick definition of a determinant,
but it has the downside that it requires showing that there is exactly one
determinant function (which is often not easier than our approach).}; this is
a valid approach, but I prefer to use Definition \ref{def.det} instead, since
it is more explicit.

\begin{exercise}
\label{exe.ps4.6k}Let $n\in\mathbb{N}$. Let $A$ be an $n\times n$-matrix.
Prove the following:

\textbf{(a)} If we add a scalar multiple of a row of $A$ to another row of
$A$, then the determinant of $A$ does not change. (A \textit{scalar multiple}
of a row vector $v$ means a row vector of the form $\lambda v$, where
$\lambda\in\mathbb{K}$.)

\textbf{(b)} If we add a scalar multiple of a column of $A$ to another column
of $A$, then the determinant of $A$ does not change. (A \textit{scalar
multiple} of a column vector $v$ means a column vector of the form $\lambda
v$, where $\lambda\in\mathbb{K}$.)
\end{exercise}

\begin{example}
Let us visualize Exercise \ref{exe.ps4.6k} \textbf{(a)}. Set $n=3$ and
$A=\left(
\begin{array}
[c]{ccc}%
a & b & c\\
d & e & f\\
g & h & i
\end{array}
\right)  $. If we add $-2$ times the second row of $A$ to the first row of
$A$, then we obtain the matrix $\left(
\begin{array}
[c]{ccc}%
a+\left(  -2\right)  d & b+\left(  -2\right)  e & c+\left(  -2\right)  f\\
d & e & f\\
g & h & i
\end{array}
\right)  $. Exercise \ref{exe.ps4.6k} \textbf{(a)} now claims that this new
matrix has the same determinant as $A$ (because $-2$ times the second row of
$A$ is a scalar multiple of the second row of $A$).

Notice the word \textquotedblleft another\textquotedblright\ in Exercise
\ref{exe.ps4.6k}. Adding a scalar multiple of a row of $A$ to \textbf{the
same} row of $A$ will likely change the determinant.
\end{example}

\Needspace{30\baselineskip}

\begin{remark}
\label{exam.xi+yj.2}Exercise \ref{exe.ps4.6k} lets us prove the claim of
Example \ref{exam.xi+yj} in a much simpler way.

Namely, let $n$ and $x_{1},x_{2},\ldots,x_{n}$ and $y_{1},y_{2},\ldots,y_{n}$
and $A$ be as in Example \ref{exam.xi+yj}. Assume that $n\geq3$. We want to
show that $\det A=0$.

The matrix $A$ has at least three rows (since $n\geq3$), and looks as follows:%
\[
A=\left(
\begin{array}
[c]{cccccc}%
x_{1}+y_{1} & x_{1}+y_{2} & x_{1}+y_{3} & x_{1}+y_{4} & \cdots & x_{1}+y_{n}\\
x_{2}+y_{1} & x_{2}+y_{2} & x_{2}+y_{3} & x_{2}+y_{4} & \cdots & x_{2}+y_{n}\\
x_{3}+y_{1} & x_{3}+y_{2} & x_{3}+y_{3} & x_{3}+y_{4} & \cdots & x_{3}+y_{n}\\
x_{4}+y_{1} & x_{4}+y_{2} & x_{4}+y_{3} & x_{4}+y_{4} & \cdots & x_{4}+y_{n}\\
\vdots & \vdots & \vdots & \vdots & \ddots & \vdots\\
x_{n}+y_{1} & x_{n}+y_{2} & x_{n}+y_{3} & x_{n}+y_{4} & \cdots & x_{n}+y_{n}%
\end{array}
\right)
\]
(where the presence of terms like $x_{4}$ and $y_{4}$ does not mean that the
variables $x_{4}$ and $y_{4}$ exist, in the same way as one can write
\textquotedblleft$x_{1},x_{2},\ldots,x_{k}$\textquotedblright\ even if $k=1$
or $k=0$). Thus, if we subtract the first row of $A$ from the second row of
$A$, then we obtain the matrix%
\[
A^{\prime}=\left(
\begin{array}
[c]{cccccc}%
x_{1}+y_{1} & x_{1}+y_{2} & x_{1}+y_{3} & x_{1}+y_{4} & \cdots & x_{1}+y_{n}\\
x_{2}-x_{1} & x_{2}-x_{1} & x_{2}-x_{1} & x_{2}-x_{1} & \cdots & x_{2}-x_{1}\\
x_{3}+y_{1} & x_{3}+y_{2} & x_{3}+y_{3} & x_{3}+y_{4} & \cdots & x_{3}+y_{n}\\
x_{4}+y_{1} & x_{4}+y_{2} & x_{4}+y_{3} & x_{4}+y_{4} & \cdots & x_{4}+y_{n}\\
\vdots & \vdots & \vdots & \vdots & \ddots & \vdots\\
x_{n}+y_{1} & x_{n}+y_{2} & x_{n}+y_{3} & x_{n}+y_{4} & \cdots & x_{n}+y_{n}%
\end{array}
\right)
\]
(because $\left(  x_{2}+y_{j}\right)  -\left(  x_{1}+y_{j}\right)
=x_{2}-x_{1}$ for every $j$). The transformation we just did (subtracting a
row from another row) does not change the determinant of the matrix (by
Exercise \ref{exe.ps4.6k} \textbf{(a)}, because subtracting a row from another
row is tantamount to adding the $\left(  -1\right)  $-multiple of the former
row to the latter), and thus we have $\det A^{\prime}=\det A$.

We notice that each entry of the second row of $A^{\prime}$ equals
$x_{2}-x_{1}$.

Next, we subtract the first row of $A^{\prime}$ from the third row of
$A^{\prime}$, and obtain the matrix%
\[
A^{\prime\prime}=\left(
\begin{array}
[c]{cccccc}%
x_{1}+y_{1} & x_{1}+y_{2} & x_{1}+y_{3} & x_{1}+y_{4} & \cdots & x_{1}+y_{n}\\
x_{2}-x_{1} & x_{2}-x_{1} & x_{2}-x_{1} & x_{2}-x_{1} & \cdots & x_{2}-x_{1}\\
x_{3}-x_{1} & x_{3}-x_{1} & x_{3}-x_{1} & x_{3}-x_{1} & \cdots & x_{3}-x_{1}\\
x_{4}+y_{1} & x_{4}+y_{2} & x_{4}+y_{3} & x_{4}+y_{4} & \cdots & x_{4}+y_{n}\\
\vdots & \vdots & \vdots & \vdots & \ddots & \vdots\\
x_{n}+y_{1} & x_{n}+y_{2} & x_{n}+y_{3} & x_{n}+y_{4} & \cdots & x_{n}+y_{n}%
\end{array}
\right)  .
\]
Again, the determinant is unchanged (because of Exercise \ref{exe.ps4.6k}
\textbf{(a)}), so we have $\det A^{\prime\prime}=\det A^{\prime}=\det A$.

We notice that each entry of the second row of $A^{\prime\prime}$ equals
$x_{2}-x_{1}$ (indeed, these entries have been copied over from $A^{\prime}$),
and that each entry of the third row of $A^{\prime\prime}$ equals $x_{3}%
-x_{1}$.

Next, we subtract the first column of $A^{\prime\prime}$ from each of the
other columns of $A^{\prime\prime}$. This gives us the matrix%
\begin{equation}
A^{\prime\prime\prime}=\left(
\begin{array}
[c]{cccccc}%
x_{1}+y_{1} & y_{2}-y_{1} & y_{3}-y_{1} & y_{4}-y_{1} & \cdots & y_{n}-y_{1}\\
x_{2}-x_{1} & 0 & 0 & 0 & \cdots & 0\\
x_{3}-x_{1} & 0 & 0 & 0 & \cdots & 0\\
x_{4}+y_{1} & y_{2}-y_{1} & y_{3}-y_{1} & y_{4}-y_{1} & \cdots & y_{n}-y_{1}\\
\vdots & \vdots & \vdots & \vdots & \ddots & \vdots\\
x_{n}+y_{1} & y_{2}-y_{1} & y_{3}-y_{1} & y_{4}-y_{1} & \cdots & y_{n}-y_{1}%
\end{array}
\right)  . \label{eq.exam.xi+yj.2.4}%
\end{equation}
This step, again, has not changed the determinant (because Exercise
\ref{exe.ps4.6k} \textbf{(b)} shows that subtracting a column from another
column does not change the determinant, and what we did was doing $n-1$ such
transformations). Thus, $\det A^{\prime\prime\prime}=\det A^{\prime\prime
}=\det A$.

Now, let us write the matrix $A^{\prime\prime\prime}$ in the form
$A^{\prime\prime\prime}=\left(  a_{i,j}^{\prime\prime\prime}\right)  _{1\leq
i\leq n,\ 1\leq j\leq n}$. (Thus, $a_{i,j}^{\prime\prime\prime}$ is the
$\left(  i,j\right)  $-th entry of $A^{\prime\prime\prime}$ for every $\left(
i,j\right)  $.) Then, (\ref{eq.det.eq.1}) (applied to $A^{\prime\prime\prime}$
instead of $A$) yields%
\begin{equation}
\det A^{\prime\prime\prime}=\sum_{\sigma\in S_{n}}\left(  -1\right)  ^{\sigma
}a_{1,\sigma\left(  1\right)  }^{\prime\prime\prime}a_{2,\sigma\left(
2\right)  }^{\prime\prime\prime}\cdots a_{n,\sigma\left(  n\right)  }%
^{\prime\prime\prime}. \label{eq.exam.xi+yj.2.5}%
\end{equation}
I claim that
\begin{equation}
a_{1,\sigma\left(  1\right)  }^{\prime\prime\prime}a_{2,\sigma\left(
2\right)  }^{\prime\prime\prime}\cdots a_{n,\sigma\left(  n\right)  }%
^{\prime\prime\prime}=0\ \ \ \ \ \ \ \ \ \ \text{for every }\sigma\in S_{n}.
\label{eq.exam.xi+yj.2.6}%
\end{equation}

[\textit{Proof of (\ref{eq.exam.xi+yj.2.6}):} Let $\sigma\in S_{n}$. Then,
$\sigma$ is injective, and thus $\sigma\left(  2\right)  \neq\sigma\left(
3\right)  $. Therefore, at least one of the integers $\sigma\left(  2\right)
$ and $\sigma\left(  3\right)  $ must be $\neq1$ (because otherwise, we would
have $\sigma\left(  2\right)  =1=\sigma\left(  3\right)  $, contradicting
$\sigma\left(  2\right)  \neq\sigma\left(  3\right)  $). We WLOG assume that
$\sigma\left(  2\right)  \neq1$. But a look at (\ref{eq.exam.xi+yj.2.4})
reveals that all entries of the second row of $A^{\prime\prime\prime}$ are
zero except for the first entry. Thus, $a_{2,j}^{\prime\prime\prime}=0$ for
every $j\neq1$. Applied to $j=\sigma\left(  2\right)  $, this yields
$a_{2,\sigma\left(  2\right)  }^{\prime\prime\prime}=0$ (since $\sigma\left(
2\right)  \neq1$). Hence, $a_{1,\sigma\left(  1\right)  }^{\prime\prime\prime
}a_{2,\sigma\left(  2\right)  }^{\prime\prime\prime}\cdots a_{n,\sigma\left(
n\right)  }^{\prime\prime\prime}=0$ (because if $0$ appears as a factor in a
product, then the whole product must be $0$). This proves
(\ref{eq.exam.xi+yj.2.6}).]

Now, (\ref{eq.exam.xi+yj.2.5}) becomes%
\[
\det A^{\prime\prime\prime}=\sum_{\sigma\in S_{n}}\left(  -1\right)  ^{\sigma
}\underbrace{a_{1,\sigma\left(  1\right)  }^{\prime\prime\prime}%
a_{2,\sigma\left(  2\right)  }^{\prime\prime\prime}\cdots a_{n,\sigma\left(
n\right)  }^{\prime\prime\prime}}_{\substack{=0\\\text{(by
(\ref{eq.exam.xi+yj.2.6}))}}}=\sum_{\sigma\in S_{n}}\left(  -1\right)
^{\sigma}0=0.
\]
Compared with $\det A^{\prime\prime\prime}=\det A$, this yields $\det A=0$.
Thus, $\det A=0$ is proven again.
\end{remark}

\begin{remark}
\label{exam.xmax}Here is another example for the use of Exercise
\ref{exe.ps4.6k}.

Let $n\in\mathbb{N}$. Let $x_{1},x_{2},\ldots,x_{n}$ be $n$ elements of
$\mathbb{K}$. Let $A$ be the matrix $\left(  x_{\max\left\{  i,j\right\}
}\right)  _{1\leq i\leq n,\ 1\leq j\leq n}$. (Recall that $\max S$ denotes the
greatest element of a nonempty set $S$.)

For example, if $n=4$, then%
\[
A=\left(
\begin{array}
[c]{cccc}%
x_{1} & x_{2} & x_{3} & x_{4}\\
x_{2} & x_{2} & x_{3} & x_{4}\\
x_{3} & x_{3} & x_{3} & x_{4}\\
x_{4} & x_{4} & x_{4} & x_{4}%
\end{array}
\right)  .
\]

We want to find $\det A$. First, let us subtract the first row of $A$ from
each of the other rows of $A$. Thus we obtain a new matrix $A^{\prime}$. The
determinant has not changed (according to Exercise \ref{exe.ps4.6k}
\textbf{(a)}); i.e., we have $\det A^{\prime}=\det A$. Here is how $A^{\prime
}$ looks like in the case when $n=4$:%
\begin{equation}
A^{\prime}=\left(
\begin{array}
[c]{cccc}%
x_{1} & x_{2} & x_{3} & x_{4}\\
x_{2}-x_{1} & 0 & 0 & 0\\
x_{3}-x_{1} & x_{3}-x_{2} & 0 & 0\\
x_{4}-x_{1} & x_{4}-x_{2} & x_{4}-x_{3} & 0
\end{array}
\right)  . \label{eq.exam.xmax.5}%
\end{equation}
Notice the many zeroes; zeroes are useful when computing determinants. To
generalize the pattern we see on (\ref{eq.exam.xmax.5}), we write the matrix
$A^{\prime}$ in the form $A^{\prime}=\left(  a_{i,j}^{\prime}\right)  _{1\leq
i\leq n,\ 1\leq j\leq n}$ (so that $a_{i,j}^{\prime}$ is the $\left(
i,j\right)  $-th entry of $A^{\prime}$ for every $\left(  i,j\right)  $).
Then, for every $\left(  i,j\right)  \in\left\{  1,2,\ldots,n\right\}  ^{2}$,
we have%
\begin{equation}
a_{i,j}^{\prime}=
\begin{cases}
x_{\max\left\{  i,j\right\}  }, & \text{if }i=1;\\
x_{\max\left\{  i,j\right\}  }-x_{\max\left\{  1,j\right\}  }, & \text{if }i>1
\end{cases}
\label{eq.exam.xmax.4}%
\end{equation}
(since we obtained the matrix $A^{\prime}$ by subtracting the first row of $A$
from each of the other rows of $A$). Hence, for every $\left(  i,j\right)
\in\left\{  1,2,\ldots,n\right\}  ^{2}$ satisfying $1<i\leq j$, we have%
\begin{align}
a_{i,j}^{\prime}  &  =x_{\max\left\{  i,j\right\}  }-x_{\max\left\{
1,j\right\}  }=x_{j}-x_{j}\ \ \ \ \ \ \ \ \ \ \left(
\begin{array}
[c]{c}%
\text{since }\max\left\{  i,j\right\}  =j\text{ (because }i\leq j\text{)}\\
\text{and }\max\left\{  1,j\right\}  =j\text{ (because }1<j\text{)}%
\end{array}
\right) \nonumber\\
&  =0. \label{eq.exam.xmax.6}%
\end{align}
This is the general explanation for the six $0$'s in (\ref{eq.exam.xmax.5}).
We notice also that the first row of the matrix $A^{\prime}$ is $\left(
x_{1},x_{2},\ldots,x_{n}\right)  $.

Now, we want to transform $A^{\prime}$ further. Namely, we first swap the
first row with the second row; then we swap the second row (which used to be
the first row) with the third row; then, the third row with the fourth row,
and so on, until we finally swap the $\left(  n-1\right)  $-th row with the
$n$-th row. As a result of these $n-1$ swaps, the first row has moved all the
way down to the bottom, past all the other rows. We denote the resulting
matrix by $A^{\prime\prime}$. For instance, if $n=4$, then%
\begin{equation}
A^{\prime\prime}=\left(
\begin{array}
[c]{cccc}%
x_{2}-x_{1} & 0 & 0 & 0\\
x_{3}-x_{1} & x_{3}-x_{2} & 0 & 0\\
x_{4}-x_{1} & x_{4}-x_{2} & x_{4}-x_{3} & 0\\
x_{1} & x_{2} & x_{3} & x_{4}%
\end{array}
\right)  . \label{eq.exam.xmax.9}%
\end{equation}
This is a lower-triangular matrix. To see that this holds in the general case,
we write the matrix $A^{\prime\prime}$ in the form $A^{\prime\prime}=\left(
a_{i,j}^{\prime\prime}\right)  _{1\leq i\leq n,\ 1\leq j\leq n}$ (so that
$a_{i,j}^{\prime\prime}$ is the $\left(  i,j\right)  $-th entry of
$A^{\prime\prime}$ for every $\left(  i,j\right)  $). Then, for every $\left(
i,j\right)  \in\left\{  1,2,\ldots,n\right\}  ^{2}$, we have%
\begin{equation}
a_{i,j}^{\prime\prime}=
\begin{cases}
a_{i+1,j}^{\prime}, & \text{if }i<n;\\
a_{1,j}^{\prime}, & \text{if }i=n
\end{cases}
\label{eq.exam.xmax.10}%
\end{equation}
(because the first row of $A^{\prime}$ has become the $n$-th row of
$A^{\prime\prime}$, whereas every other row has moved up one step). In
particular, for every $\left(  i,j\right)  \in\left\{  1,2,\ldots,n\right\}
^{2}$ satisfying $1\leq i<j\leq n$, we have%
\begin{align*}
a_{i,j}^{\prime\prime}  &  =
\begin{cases}
a_{i+1,j}^{\prime}, & \text{if }i<n;\\
a_{1,j}^{\prime}, & \text{if }i=n
\end{cases}
=a_{i+1,j}^{\prime}\ \ \ \ \ \ \ \ \ \ \left(  \text{since }i<j\leq n\right)
\\
&  =0\ \ \ \ \ \ \ \ \ \ \left(
\begin{array}
[c]{c}%
\text{by (\ref{eq.exam.xmax.6}), applied to }i+1\text{ instead of }i\\
\text{(because }i<j\text{ yields }i+1\leq j\text{)}%
\end{array}
\right)  .
\end{align*}
This shows that $A^{\prime\prime}$ is indeed lower-triangular. Hence, Exercise
\ref{exe.ps4.3} (applied to $A^{\prime\prime}$ and $a_{i,j}^{\prime\prime}$
instead of $A$ and $a_{i,j}$) shows that $\det A^{\prime\prime}=a_{1,1}%
^{\prime\prime}a_{2,2}^{\prime\prime}\cdots a_{n,n}^{\prime\prime}$.

Using (\ref{eq.exam.xmax.10}) and (\ref{eq.exam.xmax.4}), it is easy to see
that every $i\in\left\{  1,2,\ldots,n\right\}  $ satisfies%
\begin{equation}
a_{i,i}^{\prime\prime}=
\begin{cases}
x_{i+1}-x_{i}, & \text{if }i<n;\\
x_{n}, & \text{if }i=n
\end{cases}
. \label{eq.exam.xmax.15}%
\end{equation}
(This is precisely the pattern you would guess from the diagonal entries in
(\ref{eq.exam.xmax.9}).) Now, multiplying the equalities
(\ref{eq.exam.xmax.15}) for all $i\in\left\{  1,2,\ldots,n\right\}  $, we
obtain $a_{1,1}^{\prime\prime}a_{2,2}^{\prime\prime}\cdots a_{n,n}%
^{\prime\prime}=\left(  x_{2}-x_{1}\right)  \left(  x_{3}-x_{2}\right)
\cdots\left(  x_{n}-x_{n-1}\right)  x_{n}$. Thus,%
\begin{equation}
\det A^{\prime\prime}=a_{1,1}^{\prime\prime}a_{2,2}^{\prime\prime}\cdots
a_{n,n}^{\prime\prime}=\left(  x_{2}-x_{1}\right)  \left(  x_{3}-x_{2}\right)
\cdots\left(  x_{n}-x_{n-1}\right)  x_{n}. \label{eq.exam.xmax.17}%
\end{equation}

But we want $\det A$, not $\det A^{\prime\prime}$. First, let us find $\det
A^{\prime}$. Recall that $A^{\prime\prime}$ was obtained from $A^{\prime}$ by
swapping rows, repeatedly -- namely, $n-1$ times. Every time we swap two rows
in a matrix, its determinant gets multiplied by $-1$ (because of Exercise
\ref{exe.ps4.6} \textbf{(a)}). Hence, $n-1$ such swaps cause the determinant
to be multiplied by $\left(  -1\right)  ^{n-1}$. Since $A^{\prime\prime}$ was
obtained from $A^{\prime}$ by $n-1$ such swaps, we thus conclude that $\det
A^{\prime\prime}=\left(  -1\right)  ^{n-1}\det A^{\prime}$, so that%
\begin{align*}
\det A^{\prime}  &  =\underbrace{\dfrac{1}{\left(  -1\right)  ^{n-1}}%
}_{=\left(  -1\right)  ^{n-1}}\underbrace{\det A^{\prime\prime}}_{=\left(
x_{2}-x_{1}\right)  \left(  x_{3}-x_{2}\right)  \cdots\left(  x_{n}%
-x_{n-1}\right)  x_{n}}\\
&  =\left(  -1\right)  ^{n-1}\left(  x_{2}-x_{1}\right)  \left(  x_{3}%
-x_{2}\right)  \cdots\left(  x_{n}-x_{n-1}\right)  x_{n}.
\end{align*}

Finally, recall that $\det A^{\prime}=\det A$, so that%
\[
\det A=\det A^{\prime}=\left(  -1\right)  ^{n-1}\left(  x_{2}-x_{1}\right)
\left(  x_{3}-x_{2}\right)  \cdots\left(  x_{n}-x_{n-1}\right)  x_{n}.
\]

\end{remark}

\subsection{$\det\left(  AB\right)  $}

Next, a lemma that will come handy in a more important proof:

\begin{lemma}
\label{lem.det.sigma}Let $n\in\mathbb{N}$. Let $\left[  n\right]  $ denote the
set $\left\{  1,2,\ldots,n\right\}  $. Let $\kappa:\left[  n\right]
\rightarrow\left[  n\right]  $ be a map. Let $B=\left(  b_{i,j}\right)
_{1\leq i\leq n,\ 1\leq j\leq n}$ be an $n\times n$-matrix. Let $B_{\kappa}$
be the $n\times n$-matrix $\left(  b_{\kappa\left(  i\right)  ,j}\right)
_{1\leq i\leq n,\ 1\leq j\leq n}$.

\textbf{(a)} If $\kappa\in S_{n}$, then $\det\left(  B_{\kappa}\right)
=\left(  -1\right)  ^{\kappa}\cdot\det B$.

\textbf{(b)} If $\kappa\notin S_{n}$, then $\det\left(  B_{\kappa}\right)  =0$.
\end{lemma}

\begin{remark}
Lemma \ref{lem.det.sigma} \textbf{(a)} simply says that if we permute the rows
of a square matrix, then its determinant gets multiplied by the sign of the
permutation used. For instance, let $n=3$ and $B=\left(
\begin{array}
[c]{ccc}%
a & b & c\\
d & e & f\\
g & h & i
\end{array}
\right)  $. If $\kappa$ is the permutation $\left(  2,3,1\right)  $ (in
one-line notation), then $B_{\kappa}=\left(
\begin{array}
[c]{ccc}%
d & e & f\\
g & h & i\\
a & b & c
\end{array}
\right)  $, and Lemma \ref{lem.det.sigma} \textbf{(a)} says that $\det\left(
B_{\kappa}\right)  =\underbrace{\left(  -1\right)  ^{\kappa}}_{=1}\cdot\det
B=\det B$.

Of course, a similar result holds for permutations of columns.
\end{remark}

\begin{remark}
Exercise \ref{exe.ps4.6} \textbf{(a)} is a particular case of Lemma
\ref{lem.det.sigma} \textbf{(a)}. Indeed, if $B$ is an $n\times n$-matrix
obtained from $A$ by swapping the $u$-th and the $v$-th row (where $u$ and $v$
are two distinct elements of $\left\{  1,2,\ldots,n\right\}  $), then
$B=\left(  a_{t_{u,v}\left(  i\right)  ,j}\right)  _{1\leq i\leq n,\ 1\leq
j\leq n}$ (where $A$ is written in the form $A=\left(  a_{i,j}\right)  _{1\leq
i\leq n,\ 1\leq j\leq n}$).
\end{remark}

\begin{proof}
[Proof of Lemma \ref{lem.det.sigma}.]Recall that $S_{n}$ is the set of all
permutations of $\left\{  1,2,\ldots,n\right\}  $. In other words, $S_{n}$ is
the set of all permutations of $\left[  n\right]  $ (since $\left[  n\right]
=\left\{  1,2,\ldots,n\right\}  $). In other words, $S_{n}$ is the set of all
bijective maps $\left[  n\right]  \rightarrow\left[  n\right]  $.

\textbf{(a)} Assume that $\kappa\in S_{n}$. We define a map $\Phi
:S_{n}\rightarrow S_{n}$ by%
\[
\Phi\left(  \sigma\right)  =\sigma\circ\kappa\ \ \ \ \ \ \ \ \ \ \text{for
every }\sigma\in S_{n}.
\]
We also define a map $\Psi:S_{n}\rightarrow S_{n}$ by%
\[
\Psi\left(  \sigma\right)  =\sigma\circ\kappa^{-1}%
\ \ \ \ \ \ \ \ \ \ \text{for every }\sigma\in S_{n}.
\]
The maps $\Phi$ and $\Psi$ are mutually inverse\footnote{\textit{Proof.} Every
$\sigma\in S_{n}$ satisfies%
\begin{align*}
\left(  \Psi\circ\Phi\right)  \left(  \sigma\right)   &  =\Psi\left(
\underbrace{\Phi\left(  \sigma\right)  }_{=\sigma\circ\kappa}\right)
=\Psi\left(  \sigma\circ\kappa\right)  =\sigma\circ\underbrace{\kappa
\circ\kappa^{-1}}_{=\operatorname*{id}}\ \ \ \ \ \ \ \ \ \ \left(  \text{by
the definition of }\Psi\right) \\
&  =\sigma=\operatorname*{id}\left(  \sigma\right)  .
\end{align*}
Thus, $\Psi\circ\Phi=\operatorname*{id}$. Similarly, $\Phi\circ\Psi
=\operatorname*{id}$. Combined with $\Psi\circ\Phi=\operatorname*{id}$, this
yields that the maps $\Phi$ and $\Psi$ are mutually inverse, qed.}. Hence, the
map $\Phi$ is a bijection.

We have $B=\left(  b_{i,j}\right)  _{1\leq i\leq n,\ 1\leq j\leq n}$. Hence,
(\ref{eq.det.eq.2}) (applied to $B$ and $b_{i,j}$ instead of $A$ and $a_{i,j}%
$) yields%
\begin{equation}
\det B=\sum_{\sigma\in S_{n}}\left(  -1\right)  ^{\sigma}\underbrace{\prod
_{i=1}^{n}}_{=\prod_{i\in\left[  n\right]  }}b_{i,\sigma\left(  i\right)
}=\sum_{\sigma\in S_{n}}\left(  -1\right)  ^{\sigma}\prod_{i\in\left[
n\right]  }b_{i,\sigma\left(  i\right)  }. \label{pf.lem.det.sigma.a.detB}%
\end{equation}

Now, $B_{\kappa}=\left(  b_{\kappa\left(  i\right)  ,j}\right)  _{1\leq i\leq
n,\ 1\leq j\leq n}$. Hence, (\ref{eq.det.eq.2}) (applied to $B_{\kappa}$ and
$b_{\kappa\left(  i\right)  ,j}$ instead of $A$ and $a_{i,j}$) yields%
\begin{align}
\det\left(  B_{\kappa}\right)   &  =\sum_{\sigma\in S_{n}}\left(  -1\right)
^{\sigma}\underbrace{\prod_{i=1}^{n}}_{=\prod_{i\in\left[  n\right]  }%
}b_{\kappa\left(  i\right)  ,\sigma\left(  i\right)  }=\sum_{\sigma\in S_{n}%
}\left(  -1\right)  ^{\sigma}\prod_{i\in\left[  n\right]  }b_{\kappa\left(
i\right)  ,\sigma\left(  i\right)  }\nonumber\\
&  =\sum_{\sigma\in S_{n}}\left(  -1\right)  ^{\Phi\left(  \sigma\right)
}\prod_{i\in\left[  n\right]  }b_{\kappa\left(  i\right)  ,\left(  \Phi\left(
\sigma\right)  \right)  \left(  i\right)  } \label{pf.lem.det.sigma.a.1}%
\end{align}
(here, we have substituted $\Phi\left(  \sigma\right)  $ for $\sigma$ in the
sum, since $\Phi$ is a bijection).

But every $\sigma\in S_{n}$ satisfies $\left(  -1\right)  ^{\Phi\left(
\sigma\right)  }=\left(  -1\right)  ^{\kappa}\cdot\left(  -1\right)  ^{\sigma
}$\ \ \ \ \footnote{\textit{Proof.} Let $\sigma\in S_{n}$. Then, $\Phi\left(
\sigma\right)  =\sigma\circ\kappa$, so that%
\begin{align*}
\left(  -1\right)  ^{\Phi\left(  \sigma\right)  }  &  =\left(  -1\right)
^{\sigma\circ\kappa}=\left(  -1\right)  ^{\sigma}\cdot\left(  -1\right)
^{\kappa}\ \ \ \ \ \ \ \ \ \ \left(  \text{by (\ref{eq.sign.prod}), applied to
}\tau=\kappa\right) \\
&  =\left(  -1\right)  ^{\kappa}\cdot\left(  -1\right)  ^{\sigma},
\end{align*}
qed.} and $\prod_{i\in\left[  n\right]  }b_{\kappa\left(  i\right)  ,\left(
\Phi\left(  \sigma\right)  \right)  \left(  i\right)  }=\prod_{i\in\left[
n\right]  }b_{i,\sigma\left(  i\right)  }$\ \ \ \ \footnote{\textit{Proof.}
Let $\sigma\in S_{n}$. We have $\Phi\left(  \sigma\right)  =\sigma\circ\kappa
$. Thus, for every $i\in\left[  n\right]  $, we have $\left(  \Phi\left(
\sigma\right)  \right)  \left(  i\right)  =\left(  \sigma\circ\kappa\right)
\left(  i\right)  =\sigma\left(  \kappa\left(  i\right)  \right)  $. Hence,
$\prod_{i\in\left[  n\right]  }b_{\kappa\left(  i\right)  ,\left(  \Phi\left(
\sigma\right)  \right)  \left(  i\right)  }=\prod_{i\in\left[  n\right]
}b_{\kappa\left(  i\right)  ,\sigma\left(  \kappa\left(  i\right)  \right)  }%
$.
\par
But $\kappa\in S_{n}$. In other words, $\kappa$ is a permutation of the set
$\left\{  1,2,\ldots,n\right\}  =\left[  n\right]  $, hence a bijection from
$\left[  n\right]  $ to $\left[  n\right]  $. Therefore, we can substitute
$\kappa\left(  i\right)  $ for $i$ in the product $\prod_{i\in\left[
n\right]  }b_{i,\sigma\left(  i\right)  }$. We thus obtain $\prod_{i\in\left[
n\right]  }b_{i,\sigma\left(  i\right)  }=\prod_{i\in\left[  n\right]
}b_{\kappa\left(  i\right)  ,\sigma\left(  \kappa\left(  i\right)  \right)  }%
$. Comparing this with $\prod_{i\in\left[  n\right]  }b_{\kappa\left(
i\right)  ,\left(  \Phi\left(  \sigma\right)  \right)  \left(  i\right)
}=\prod_{i\in\left[  n\right]  }b_{\kappa\left(  i\right)  ,\sigma\left(
\kappa\left(  i\right)  \right)  }$, we obtain $\prod_{i\in\left[  n\right]
}b_{\kappa\left(  i\right)  ,\left(  \Phi\left(  \sigma\right)  \right)
\left(  i\right)  }=\prod_{i\in\left[  n\right]  }b_{i,\sigma\left(  i\right)
}$, qed.}. Thus, (\ref{pf.lem.det.sigma.a.1}) becomes%
\begin{align*}
\det\left(  B_{\kappa}\right)   &  =\sum_{\sigma\in S_{n}}\underbrace{\left(
-1\right)  ^{\Phi\left(  \sigma\right)  }}_{=\left(  -1\right)  ^{\kappa}%
\cdot\left(  -1\right)  ^{\sigma}}\underbrace{\prod_{i\in\left[  n\right]
}b_{\kappa\left(  i\right)  ,\left(  \Phi\left(  \sigma\right)  \right)
\left(  i\right)  }}_{=\prod_{i\in\left[  n\right]  }b_{i,\sigma\left(
i\right)  }}=\sum_{\sigma\in S_{n}}\left(  -1\right)  ^{\kappa}\cdot\left(
-1\right)  ^{\sigma}\prod_{i\in\left[  n\right]  }b_{i,\sigma\left(  i\right)
}\\
&  =\left(  -1\right)  ^{\kappa}\cdot\underbrace{\sum_{\sigma\in S_{n}}\left(
-1\right)  ^{\sigma}\prod_{i\in\left[  n\right]  }b_{i,\sigma\left(  i\right)
}}_{\substack{=\det B\\\text{(by (\ref{pf.lem.det.sigma.a.detB}))}}}=\left(
-1\right)  ^{\kappa}\cdot\det B.
\end{align*}
This proves Lemma \ref{lem.det.sigma} \textbf{(a)}.

\textbf{(b)} Assume that $\kappa\notin S_{n}$.

\begin{vershort}
The following fact is well-known: If $U$ is a finite set, then every injective
map $U\rightarrow U$ is bijective\footnote{\textit{Proof.} Let $U$ be a finite
set, and let $f$ be an injective map $U\rightarrow U$. We must show that $f$
is bijective.
\par
Since $f$ is injective, we have $\left\vert f\left(  U\right)  \right\vert
=\left\vert U\right\vert $. Thus, $f\left(  U\right)  $ is a subset of $U$
which has size $\left\vert U\right\vert $. But the only such subset is $U$
itself (since $U$ is a finite set). Therefore, $f\left(  U\right)  $ must be
$U$ itself. In other words, the map $f$ is surjective. Hence, $f$ is bijective
(since $f$ is injective and surjective), qed.}. We can apply this to
$U=\left[  n\right]  $, and thus conclude that every injective map $\left[
n\right]  \rightarrow\left[  n\right]  $ is bijective. Therefore, if the map
$\kappa:\left[  n\right]  \rightarrow\left[  n\right]  $ were injective, then
$\kappa$ would be bijective and therefore would be an element of $S_{n}$
(since $S_{n}$ is the set of all bijective maps $\left[  n\right]
\rightarrow\left[  n\right]  $); but this would contradict the fact that
$\kappa\notin S_{n}$. Hence, the map $\kappa:\left[  n\right]  \rightarrow
\left[  n\right]  $ cannot be injective. Therefore, there exist two distinct
elements $a$ and $b$ of $\left[  n\right]  $ such that $\kappa\left(
a\right)  =\kappa\left(  b\right)  $. Consider these $a$ and $b$.
\end{vershort}

\begin{verlong}
Assume (for the sake of contradiction) that the map $\kappa:\left[  n\right]
\rightarrow\left[  n\right]  $ is injective. The set $\left[  n\right]  $ is
finite and satisfies $\left\vert \left[  n\right]  \right\vert \geq\left\vert
\left[  n\right]  \right\vert $. Hence, Lemma \ref{lem.jectivity.pigeon-inj}
(applied to $\left[  n\right]  $, $\left[  n\right]  $ and $\kappa$ instead of
$U$, $V$ and $f$) shows that we have the following logical equivalence:%
\[
\left(  \kappa\text{ is injective}\right)  \ \Longleftrightarrow\ \left(
\kappa\text{ is bijective}\right)  .
\]
Thus, $\kappa$ is bijective (since $\kappa$ is injective). In other words,
$\kappa$ is a bijection $\left[  n\right]  \rightarrow\left[  n\right]  $. In
other words, $\kappa$ is an element of $S_{n}$ (since $S_{n}$ is the set of
all bijective maps $\left[  n\right]  \rightarrow\left[  n\right]  $). But
this contradicts the fact that $\kappa\notin S_{n}$. This contradiction shows
that our assumption (that the map $\kappa:\left[  n\right]  \rightarrow\left[
n\right]  $ is injective) was wrong.

Hence, the map $\kappa:\left[  n\right]  \rightarrow\left[  n\right]  $ is not
injective. Therefore, there exist two distinct elements $a$ and $b$ of
$\left[  n\right]  $ such that $\kappa\left(  a\right)  =\kappa\left(
b\right)  $. Consider these $a$ and $b$.
\end{verlong}

Thus, $a$ and $b$ are two distinct elements of $\left[  n\right]  =\left\{
1,2,\ldots,n\right\}  $. Hence, a transposition $t_{a,b}\in S_{n}$ is defined
(see Definition \ref{def.transpos} for the definition). This transposition
satisfies $\kappa\circ t_{a,b}=\kappa$\ \ \ \ \footnote{\textit{Proof.} We are
going to show that every $i\in\left[  n\right]  $ satisfies $\left(
\kappa\circ t_{a,b}\right)  \left(  i\right)  =\kappa\left(  i\right)  $.
\par
So let $i\in\left[  n\right]  $. We shall show that $\left(  \kappa\circ
t_{a,b}\right)  \left(  i\right)  =\kappa\left(  i\right)  $.
\par
The definition of $t_{a,b}$ shows that $t_{a,b}$ is the permutation in $S_{n}$
which swaps $a$ with $b$ while leaving all other elements of $\left\{
1,2,\ldots,n\right\}  $ unchanged. In other words, we have $t_{a,b}\left(
a\right)  =b$, and $t_{a,b}\left(  b\right)  =a$, and $t_{a,b}\left(
j\right)  =j$ for every $j\in\left[  n\right]  \setminus\left\{  a,b\right\}
$.
\par
Now, we have $i\in\left[  n\right]  $. Thus, we are in one of the following
three cases:
\par
\textit{Case 1:} We have $i=a$.
\par
\textit{Case 2:} We have $i=b$.
\par
\textit{Case 3:} We have $i\in\left[  n\right]  \setminus\left\{  a,b\right\}
$.
\par
Let us first consider Case 1. In this case, we have $i=a$, so that $\left(
\kappa\circ t_{a,b}\right)  \left(  \underbrace{i}_{=a}\right)  =\left(
\kappa\circ t_{a,b}\right)  \left(  a\right)  =\kappa\left(
\underbrace{t_{a,b}\left(  a\right)  }_{=b}\right)  =\kappa\left(  b\right)
$. Compared with $\kappa\left(  \underbrace{i}_{=a}\right)  =\kappa\left(
a\right)  =\kappa\left(  b\right)  $, this yields $\left(  \kappa\circ
t_{a,b}\right)  \left(  i\right)  =\kappa\left(  i\right)  $. Thus, $\left(
\kappa\circ t_{a,b}\right)  \left(  i\right)  =\kappa\left(  i\right)  $ is
proven in Case 1.
\par
Let us next consider Case 2. In this case, we have $i=b$, so that $\left(
\kappa\circ t_{a,b}\right)  \left(  \underbrace{i}_{=b}\right)  =\left(
\kappa\circ t_{a,b}\right)  \left(  b\right)  =\kappa\left(
\underbrace{t_{a,b}\left(  b\right)  }_{=a}\right)  =\kappa\left(  a\right)
=\kappa\left(  b\right)  $. Compared with $\kappa\left(  \underbrace{i}%
_{=b}\right)  =\kappa\left(  b\right)  $, this yields $\left(  \kappa\circ
t_{a,b}\right)  \left(  i\right)  =\kappa\left(  i\right)  $. Thus, $\left(
\kappa\circ t_{a,b}\right)  \left(  i\right)  =\kappa\left(  i\right)  $ is
proven in Case 2.
\par
Let us finally consider Case 3. In this case, we have $i\in\left[  n\right]
\setminus\left\{  a,b\right\}  $. Hence, $t_{a,b}\left(  i\right)  =i$ (since
$t_{a,b}\left(  j\right)  =j$ for every $j\in\left[  n\right]  \setminus
\left\{  a,b\right\}  $). Therefore, $\left(  \kappa\circ t_{a,b}\right)
\left(  i\right)  =\kappa\left(  \underbrace{t_{a,b}\left(  i\right)  }%
_{=i}\right)  =\kappa\left(  i\right)  $. Thus, $\left(  \kappa\circ
t_{a,b}\right)  \left(  i\right)  =\kappa\left(  i\right)  $ is proven in Case
3.
\par
We now have shown $\left(  \kappa\circ t_{a,b}\right)  \left(  i\right)
=\kappa\left(  i\right)  $ in each of the three Cases 1, 2 and 3. Hence,
$\left(  \kappa\circ t_{a,b}\right)  \left(  i\right)  =\kappa\left(
i\right)  $ always holds.
\par
Now, let us forget that we fixed $i$. We thus have shown that $\left(
\kappa\circ t_{a,b}\right)  \left(  i\right)  =\kappa\left(  i\right)  $ for
every $i\in\left[  n\right]  $. In other words, $\kappa\circ t_{a,b}=\kappa$,
qed.}. Exercise \ref{exe.ps4.1ab} \textbf{(b)} (applied to $i=a$ and $j=b$)
yields $\left(  -1\right)  ^{t_{a,b}}=-1$.

Let $A_{n}$ be the set of all even permutations in $S_{n}$. Let $C_{n}$ be the
set of all odd permutations in $S_{n}$.

We have $\sigma\circ t_{a,b}\in C_{n}$ for every $\sigma\in A_{n}%
$\ \ \ \ \footnote{\textit{Proof.} Let $\sigma\in A_{n}$. Then, $\sigma$ is an
even permutation in $S_{n}$ (since $A_{n}$ is the set of all even permutations
in $S_{n}$). Hence, $\left(  -1\right)  ^{\sigma}=1$. Now, (\ref{eq.sign.prod}%
) (applied to $\tau=t_{a,b}$) yields $\left(  -1\right)  ^{\sigma\circ
t_{a,b}}=\underbrace{\left(  -1\right)  ^{\sigma}}_{=1}\cdot
\underbrace{\left(  -1\right)  ^{t_{a,b}}}_{=-1}=-1$. Thus, the permutation
$\sigma\circ t_{a,b}$ is odd. Hence, $\sigma\circ t_{a,b}$ is an odd
permutation in $S_{n}$. In other words, $\sigma\circ t_{a,b}\in C_{n}$ (since
$C_{n}$ is the set of all odd permutations in $S_{n}$), qed.}. Hence, we can
define a map $\Phi:A_{n}\rightarrow C_{n}$ by%
\[
\Phi\left(  \sigma\right)  =\sigma\circ t_{a,b}\ \ \ \ \ \ \ \ \ \ \text{for
every }\sigma\in A_{n}.
\]
Consider this map $\Phi$. Furthermore, we have $\sigma\circ\left(
t_{a,b}\right)  ^{-1}\in A_{n}$ for every $\sigma\in C_{n}$%
\ \ \ \ \footnote{\textit{Proof.} Let $\sigma\in C_{n}$. Then, $\sigma$ is an
odd permutation in $S_{n}$ (since $C_{n}$ is the set of all odd permutations
in $S_{n}$). Hence, $\left(  -1\right)  ^{\sigma}=-1$.
\par
Applying (\ref{eq.sign.inverse}) to $t_{a,b}$ instead of $\sigma$, we obtain
$\left(  -1\right)  ^{\left(  t_{a,b}\right)  ^{-1}}=\left(  -1\right)
^{t_{a,b}}=-1$. Now, (\ref{eq.sign.prod}) (applied to $\tau=\left(
t_{a,b}\right)  ^{-1}$) yields $\left(  -1\right)  ^{\sigma\circ\left(
t_{a,b}\right)  ^{-1}}=\underbrace{\left(  -1\right)  ^{\sigma}}_{=-1}%
\cdot\underbrace{\left(  -1\right)  ^{\left(  t_{a,b}\right)  ^{-1}}}%
_{=-1}=\left(  -1\right)  \cdot\left(  -1\right)  =1$. Thus, the permutation
$\sigma\circ\left(  t_{a,b}\right)  ^{-1}$ is even. Hence, $\sigma\circ\left(
t_{a,b}\right)  ^{-1}$ is an even permutation in $S_{n}$. In other words,
$\sigma\circ\left(  t_{a,b}\right)  ^{-1}\in A_{n}$ (since $A_{n}$ is the set
of all even permutations in $S_{n}$), qed.}. Thus, we can define a map
$\Psi:C_{n}\rightarrow A_{n}$ by%
\[
\Psi\left(  \sigma\right)  =\sigma\circ\left(  t_{a,b}\right)  ^{-1}%
\ \ \ \ \ \ \ \ \ \ \text{for every }\sigma\in C_{n}.
\]
Consider this map $\Psi$.

(We could have simplified our life a bit by noticing that $\left(
t_{a,b}\right)  ^{-1}=t_{a,b}$, so that the maps $\Phi$ and $\Psi$ are given
by the same formula, albeit defined on different domains. But I wanted to
demonstrate a use of (\ref{eq.sign.inverse}).)

The maps $\Phi$ and $\Psi$ are mutually inverse\footnote{\textit{Proof.} Every
$\sigma\in A_{n}$ satisfies%
\begin{align*}
\left(  \Psi\circ\Phi\right)  \left(  \sigma\right)   &  =\Psi\left(
\underbrace{\Phi\left(  \sigma\right)  }_{=\sigma\circ\tau_{a,b}}\right)
=\Psi\left(  \sigma\circ\tau_{a,b}\right)  =\sigma\circ\underbrace{\tau
_{a,b}\circ\left(  \tau_{a,b}\right)  ^{-1}}_{=\operatorname*{id}%
}\ \ \ \ \ \ \ \ \ \ \left(  \text{by the definition of }\Psi\right) \\
&  =\sigma=\operatorname*{id}\left(  \sigma\right)  .
\end{align*}
Thus, $\Psi\circ\Phi=\operatorname*{id}$. Similarly, $\Phi\circ\Psi
=\operatorname*{id}$. Combined with $\Psi\circ\Phi=\operatorname*{id}$, this
yields that the maps $\Phi$ and $\Psi$ are mutually inverse, qed.}. Hence, the
map $\Psi$ is a bijection. Moreover, every $\sigma\in C_{n}$ satisfies%
\begin{equation}
\prod_{i\in\left[  n\right]  }b_{\kappa\left(  i\right)  ,\left(  \Psi\left(
\sigma\right)  \right)  \left(  i\right)  }=\prod_{i\in\left[  n\right]
}b_{\kappa\left(  i\right)  ,\sigma\left(  i\right)  }.
\label{pf.lem.det.sigma.b.2}%
\end{equation}
\footnote{\textit{Proof of (\ref{pf.lem.det.sigma.b.2}):} Let $\sigma\in
C_{n}$. The map $t_{a,b}$ is a permutation of $\left[  n\right]  $, thus a
bijection $\left[  n\right]  \rightarrow\left[  n\right]  $. Hence, we can
substitute $t_{a,b}\left(  i\right)  $ for $i$ in the product $\prod
_{i\in\left[  n\right]  }b_{\kappa\left(  i\right)  ,\left(  \Psi\left(
\sigma\right)  \right)  \left(  i\right)  }$. Thus we obtain%
\[
\prod_{i\in\left[  n\right]  }b_{\kappa\left(  i\right)  ,\left(  \Psi\left(
\sigma\right)  \right)  \left(  i\right)  }=\prod_{i\in\left[  n\right]
}b_{\kappa\left(  t_{a,b}\left(  i\right)  \right)  ,\left(  \Psi\left(
\sigma\right)  \right)  \left(  t_{a,b}\left(  i\right)  \right)  }%
=\prod_{i\in\left[  n\right]  }b_{\kappa\left(  i\right)  ,\sigma\left(
i\right)  }%
\]
(since every $i\in\left[  n\right]  $ satisfies $\kappa\left(  t_{a,b}\left(
i\right)  \right)  =\underbrace{\left(  \kappa\circ t_{a,b}\right)  }%
_{=\kappa}\left(  i\right)  =\kappa\left(  i\right)  $ and
\[
\underbrace{\left(  \Psi\left(  \sigma\right)  \right)  }_{=\sigma\circ\left(
t_{a,b}\right)  ^{-1}}\left(  t_{a,b}\left(  i\right)  \right)  =\left(
\sigma\circ\left(  t_{a,b}\right)  ^{-1}\right)  \left(  t_{a,b}\left(
i\right)  \right)  =\sigma\left(  \underbrace{\left(  t_{a,b}\right)
^{-1}\left(  t_{a,b}\left(  i\right)  \right)  }_{=i}\right)  =\sigma\left(
i\right)
\]
). This proves (\ref{pf.lem.det.sigma.b.2}).}

We have $B_{\kappa}=\left(  b_{\kappa\left(  i\right)  ,j}\right)  _{1\leq
i\leq n,\ 1\leq j\leq n}$. Hence, (\ref{eq.det.eq.2}) (applied to $B_{\kappa}$
and $b_{\kappa\left(  i\right)  ,j}$ instead of $A$ and $a_{i,j}$) yields%
\begin{align*}
\det\left(  B_{\kappa}\right)   &  =\sum_{\sigma\in S_{n}}\left(  -1\right)
^{\sigma}\underbrace{\prod_{i=1}^{n}}_{=\prod_{i\in\left[  n\right]  }%
}b_{\kappa\left(  i\right)  ,\sigma\left(  i\right)  }=\sum_{\sigma\in S_{n}%
}\left(  -1\right)  ^{\sigma}\prod_{i\in\left[  n\right]  }b_{\kappa\left(
i\right)  ,\sigma\left(  i\right)  }\\
&  =\underbrace{\sum_{\substack{\sigma\in S_{n};\\\sigma\text{ is even}}%
}}_{\substack{=\sum_{\sigma\in A_{n}}\\\text{(since }A_{n}\text{ is
the}\\\text{set of all even}\\\text{permutations}\\\text{in }S_{n}\text{)}%
}}\underbrace{\left(  -1\right)  ^{\sigma}}_{\substack{=1\\\text{(since
}\sigma\text{ is even)}}}\prod_{i\in\left[  n\right]  }b_{\kappa\left(
i\right)  ,\sigma\left(  i\right)  }+\underbrace{\sum_{\substack{\sigma\in
S_{n};\\\sigma\text{ is odd}}}}_{\substack{=\sum_{\sigma\in C_{n}%
}\\\text{(since }C_{n}\text{ is the}\\\text{set of all odd}%
\\\text{permutations}\\\text{in }S_{n}\text{)}}}\underbrace{\left(  -1\right)
^{\sigma}}_{\substack{=-1\\\text{(since }\sigma\text{ is odd)}}}\prod
_{i\in\left[  n\right]  }b_{\kappa\left(  i\right)  ,\sigma\left(  i\right)
}\\
&  \ \ \ \ \ \ \ \ \ \ \ \ \ \ \ \ \ \ \ \ \left(
\begin{array}
[c]{c}%
\text{since every permutation }\sigma\in S_{n}\text{ is}\\
\text{either even or odd, but not both}%
\end{array}
\right) \\
&  =\sum_{\sigma\in A_{n}}\prod_{i\in\left[  n\right]  }b_{\kappa\left(
i\right)  ,\sigma\left(  i\right)  }+\sum_{\sigma\in C_{n}}\left(  -1\right)
\prod_{i\in\left[  n\right]  }b_{\kappa\left(  i\right)  ,\sigma\left(
i\right)  }\\
&  =\sum_{\sigma\in A_{n}}\prod_{i\in\left[  n\right]  }b_{\kappa\left(
i\right)  ,\sigma\left(  i\right)  }-\sum_{\sigma\in C_{n}}\prod_{i\in\left[
n\right]  }b_{\kappa\left(  i\right)  ,\sigma\left(  i\right)  }=0,
\end{align*}
since%
\begin{align*}
&  \sum_{\sigma\in A_{n}}\prod_{i\in\left[  n\right]  }b_{\kappa\left(
i\right)  ,\sigma\left(  i\right)  }\\
&  =\sum_{\sigma\in C_{n}}\underbrace{\prod_{i\in\left[  n\right]  }%
b_{\kappa\left(  i\right)  ,\left(  \Psi\left(  \sigma\right)  \right)
\left(  i\right)  }}_{\substack{=\prod_{i\in\left[  n\right]  }b_{\kappa
\left(  i\right)  ,\sigma\left(  i\right)  }\\\text{(by
(\ref{pf.lem.det.sigma.b.2}))}}}\\
&  \ \ \ \ \ \ \ \ \ \ \left(  \text{here, we have substituted }\Psi\left(
\sigma\right)  \text{ for }\sigma\text{, since the map }\Psi\text{ is a
bijection}\right) \\
&  =\sum_{\sigma\in C_{n}}\prod_{i\in\left[  n\right]  }b_{\kappa\left(
i\right)  ,\sigma\left(  i\right)  }.
\end{align*}
This proves Lemma \ref{lem.det.sigma} \textbf{(b)}.
\end{proof}

Now let us state a basic formula for products of sums in a commutative ring:

\begin{lemma}
\label{lem.prodrule}For every $n\in\mathbb{N}$, let $\left[  n\right]  $
denote the set $\left\{  1,2,\ldots,n\right\}  $.

Let $n\in\mathbb{N}$. For every $i\in\left[  n\right]  $, let $p_{i,1}%
,p_{i,2},\ldots,p_{i,m_{i}}$ be finitely many elements of $\mathbb{K}$. Then,%
\[
\prod_{i=1}^{n}\ \ \sum_{k=1}^{m_{i}}p_{i,k}=\sum_{\left(  k_{1},k_{2}%
,\ldots,k_{n}\right)  \in\left[  m_{1}\right]  \times\left[  m_{2}\right]
\times\cdots\times\left[  m_{n}\right]  }\ \ \prod_{i=1}^{n}p_{i,k_{i}}.
\]

(\textbf{Pedantic remark:} If $n=0$, then the Cartesian product $\left[
m_{1}\right]  \times\left[  m_{2}\right]  \times\cdots\times\left[
m_{n}\right]  $ has no factors; it is what is called an \textit{empty
Cartesian product}. It is understood to be a $1$-element set, and its single
element is the $0$-tuple $\left(  {}\right)  $ (also known as the empty list).)
\end{lemma}

I tend to refer to Lemma \ref{lem.prodrule} as the \textit{product rule}
(since it is related to the product rule for joint probabilities); I think it
has no really widespread name. However, it is a fundamental algebraic fact
that is used very often and tacitly (I suspect that most mathematicians have
never thought of it as being a theorem that needs to be proven). The idea
behind Lemma \ref{lem.prodrule} is that if you expand the product%
\begin{align*}
&  \prod_{i=1}^{n}\ \ \sum_{k=1}^{m_{i}}p_{i,k}\\
&  =\prod_{i=1}^{n}\left(  p_{i,1}+p_{i,2}+\cdots+p_{i,m_{i}}\right) \\
&  =\left(  p_{1,1}+p_{1,2}+\cdots+p_{1,m_{1}}\right)  \left(  p_{2,1}%
+p_{2,2}+\cdots+p_{2,m_{2}}\right)  \cdots\left(  p_{n,1}+p_{n,2}%
+\cdots+p_{n,m_{n}}\right)  ,
\end{align*}
then you get a sum of $m_{1}m_{2}\cdots m_{n}$ terms, each of which has the
form
\[
p_{1,k_{1}}p_{2,k_{2}}\cdots p_{n,k_{n}}=\prod_{i=1}^{n}p_{i,k_{i}}%
\]
for some $\left(  k_{1},k_{2},\ldots,k_{n}\right)  \in\left[  m_{1}\right]
\times\left[  m_{2}\right]  \times\cdots\times\left[  m_{n}\right]  $. (More
precisely, it is the sum of all such terms.) A formal proof of Lemma
\ref{lem.prodrule} could be obtained by induction over $n$ using the
distributivity axiom\footnote{and the observation that the $n$-tuples $\left(
k_{1},k_{2},\ldots,k_{n}\right)  \in\left[  m_{1}\right]  \times\left[
m_{2}\right]  \times\cdots\times\left[  m_{n}\right]  $ are in bijection with
the pairs $\left(  \left(  k_{1},k_{2},\ldots,k_{n-1}\right)  ,k_{n}\right)  $
of an $\left(  n-1\right)  $-tuple $\left(  k_{1},k_{2},\ldots,k_{n-1}\right)
\in\left[  m_{1}\right]  \times\left[  m_{2}\right]  \times\cdots\times\left[
m_{n-1}\right]  $ and an element $k_{n}\in\left[  m_{n}\right]  $}. For the
details (if you care about them), see the solution to the following exercise:

\begin{exercise}
\label{exe.prodrule}Prove Lemma \ref{lem.prodrule}.
\end{exercise}

\begin{remark}
\label{rmk.prodrule.ai+bi}Lemma \ref{lem.prodrule} can be regarded as a
generalization of Exercise \ref{exe.prod(ai+bi)} \textbf{(a)}. Indeed, let me
sketch how Exercise \ref{exe.prod(ai+bi)} \textbf{(a)} can be derived from
Lemma \ref{lem.prodrule}:

Let $n$, $\left(  a_{1},a_{2},\ldots,a_{n}\right)  $ and $\left(  b_{1}%
,b_{2},\ldots,b_{n}\right)  $ be as in Exercise \ref{exe.prod(ai+bi)}
\textbf{(a)}. For every $i\in\left[  n\right]  $, set $m_{i}=2$,
$p_{i,1}=a_{i}$ and $p_{i,2}=b_{i}$. Then, Lemma \ref{lem.prodrule} yields%
\begin{align}
\prod_{i=1}^{n}\left(  a_{i}+b_{i}\right)   &  =\underbrace{\sum_{\left(
k_{1},k_{2},\ldots,k_{n}\right)  \in\underbrace{\left[  2\right]
\times\left[  2\right]  \times\cdots\times\left[  2\right]  }_{n\text{
factors}}}}_{=\sum_{\left(  k_{1},k_{2},\ldots,k_{n}\right)  \in\left[
2\right]  ^{n}}}\underbrace{\prod_{i=1}^{n}p_{i,k_{i}}}_{=\left(
\prod_{\substack{i\in\left[  n\right]  ;\\k_{i}=1}}a_{i}\right)  \left(
\prod_{\substack{i\in\left[  n\right]  ;\\k_{i}=2}}b_{i}\right)  }\nonumber\\
&  =\sum_{\left(  k_{1},k_{2},\ldots,k_{n}\right)  \in\left[  2\right]  ^{n}%
}\left(  \prod_{\substack{i\in\left[  n\right]  ;\\k_{i}=1}}a_{i}\right)
\left(  \prod_{\substack{i\in\left[  n\right]  ;\\k_{i}=2}}b_{i}\right)  .
\label{eq.rmk.prodrule.ai+bi.2}%
\end{align}
But there is a bijection between the set $\left[  2\right]  ^{n}$ and the
powerset $\mathcal{P}\left(  \left[  n\right]  \right)  $ of $\left[
n\right]  $: Namely, to every $n$-tuple $\left(  k_{1},k_{2},\ldots
,k_{n}\right)  \in\left[  2\right]  ^{n}$, we can assign the set $\left\{
i\in\left[  n\right]  \ \mid\ k_{i}=1\right\}  \in\mathcal{P}\left(  \left[
n\right]  \right)  $. It is easy to see that this assignment really is a
bijection $\left[  2\right]  ^{n}\rightarrow\mathcal{P}\left(  \left[
n\right]  \right)  $, and that it furthermore has the property that every
$n$-tuple $\left(  k_{1},k_{2},\ldots,k_{n}\right)  \in\left[  2\right]  ^{n}$
satisfies%
\[
\left(  \prod_{\substack{i\in\left[  n\right]  ;\\k_{i}=1}}a_{i}\right)
\left(  \prod_{\substack{i\in\left[  n\right]  ;\\k_{i}=2}}b_{i}\right)
=\left(  \prod_{i\in I}a_{i}\right)  \left(  \prod_{i\in\left[  n\right]
\setminus I}b_{i}\right)  ,
\]
where $I$ is the image of $\left(  k_{1},k_{2},\ldots,k_{n}\right)  $ under
this bijection. Hence,%
\begin{align*}
&  \sum_{\left(  k_{1},k_{2},\ldots,k_{n}\right)  \in\left[  2\right]  ^{n}%
}\left(  \prod_{\substack{i\in\left[  n\right]  ;\\k_{i}=1}}a_{i}\right)
\left(  \prod_{\substack{i\in\left[  n\right]  ;\\k_{i}=2}}b_{i}\right) \\
&  =\sum_{I\subseteq\left[  n\right]  }\left(  \prod_{i\in I}a_{i}\right)
\left(  \prod_{i\in\left[  n\right]  \setminus I}b_{i}\right)  .
\end{align*}

Hence, (\ref{eq.rmk.prodrule.ai+bi.2}) rewrites as%
\[
\prod_{i=1}^{n}\left(  a_{i}+b_{i}\right)  =\sum_{I\subseteq\left[  n\right]
}\left(  \prod_{i\in I}a_{i}\right)  \left(  \prod_{i\in\left[  n\right]
\setminus I}b_{i}\right)  .
\]
But this is precisely the claim of Exercise \ref{exe.prod(ai+bi)} \textbf{(a)}.
\end{remark}

We shall use a corollary of Lemma \ref{lem.prodrule}:

\begin{lemma}
\label{lem.prodrule2}For every $n\in\mathbb{N}$, let $\left[  n\right]  $
denote the set $\left\{  1,2,\ldots,n\right\}  $.

Let $n\in\mathbb{N}$ and $m\in\mathbb{N}$. For every $i\in\left[  n\right]  $,
let $p_{i,1},p_{i,2},\ldots,p_{i,m}$ be $m$ elements of $\mathbb{K}$. Then,%
\[
\prod_{i=1}^{n}\ \ \sum_{k=1}^{m}p_{i,k}=\sum_{\kappa:\left[  n\right]
\rightarrow\left[  m\right]  }\ \ \prod_{i=1}^{n}p_{i,\kappa\left(  i\right)
}.
\]

\end{lemma}

\begin{proof}
[Proof of Lemma \ref{lem.prodrule2}.]For the sake of completeness, let us give
this proof.

Lemma \ref{lem.prodrule} (applied to $m_{i}=m$ for every $i\in\left[
n\right]  $) yields
\begin{equation}
\prod_{i=1}^{n}\ \ \sum_{k=1}^{m}p_{i,k}=\sum_{\left(  k_{1},k_{2}%
,\ldots,k_{n}\right)  \in\underbrace{\left[  m\right]  \times\left[  m\right]
\times\cdots\times\left[  m\right]  }_{n\text{ factors}}}\ \ \prod_{i=1}%
^{n}p_{i,k_{i}}. \label{pf.lem.prodrule2.1}%
\end{equation}

Let $\operatorname*{Map}\left(  \left[  n\right]  ,\left[  m\right]  \right)
$ denote the set of all functions from $\left[  n\right]  $ to $\left[
m\right]  $. Now, let $\Phi$ be the map from $\operatorname*{Map}\left(
\left[  n\right]  ,\left[  m\right]  \right)  $ to $\underbrace{\left[
m\right]  \times\left[  m\right]  \times\cdots\times\left[  m\right]
}_{n\text{ factors}}$ given by%
\[
\Phi\left(  \kappa\right)  =\left(  \kappa\left(  1\right)  ,\kappa\left(
2\right)  ,\ldots,\kappa\left(  n\right)  \right)
\ \ \ \ \ \ \ \ \ \ \text{for every }\kappa\in\operatorname*{Map}\left(
\left[  n\right]  ,\left[  m\right]  \right)  .
\]
So the map $\Phi$ takes a function $\kappa$ from $\left[  n\right]  $ to
$\left[  m\right]  $, and outputs the list \newline$\left(  \kappa\left(
1\right)  ,\kappa\left(  2\right)  ,\ldots,\kappa\left(  n\right)  \right)  $
of all its values. Clearly, the map $\Phi$ is injective (since a function
$\kappa\in\operatorname*{Map}\left(  \left[  n\right]  ,\left[  m\right]
\right)  $ can be reconstructed from the list $\left(  \kappa\left(  1\right)
,\kappa\left(  2\right)  ,\ldots,\kappa\left(  n\right)  \right)  =\Phi\left(
\kappa\right)  $) and surjective (since every list of $n$ elements of $\left[
m\right]  $ is the list of values of some function $\kappa\in
\operatorname*{Map}\left(  \left[  n\right]  ,\left[  m\right]  \right)  $).
Thus, $\Phi$ is bijective. Therefore, we can substitute $\Phi\left(
\kappa\right)  $ for $\left(  k_{1},k_{2},\ldots,k_{n}\right)  $ in the sum
$\sum_{\left(  k_{1},k_{2},\ldots,k_{n}\right)  \in\underbrace{\left[
m\right]  \times\left[  m\right]  \times\cdots\times\left[  m\right]
}_{n\text{ factors}}}\prod_{i=1}^{n}p_{i,k_{i}}$. In other words, we can
substitute $\left(  \kappa\left(  1\right)  ,\kappa\left(  2\right)
,\ldots,\kappa\left(  n\right)  \right)  $ for $\left(  k_{1},k_{2}%
,\ldots,k_{n}\right)  $ in this sum (since $\Phi\left(  \kappa\right)
=\left(  \kappa\left(  1\right)  ,\kappa\left(  2\right)  ,\ldots
,\kappa\left(  n\right)  \right)  $ for each $\kappa\in\operatorname*{Map}%
\left(  \left[  n\right]  ,\left[  m\right]  \right)  $). We thus obtain%
\begin{align*}
\sum_{\left(  k_{1},k_{2},\ldots,k_{n}\right)  \in\underbrace{\left[
m\right]  \times\left[  m\right]  \times\cdots\times\left[  m\right]
}_{n\text{ factors}}}\ \ \prod_{i=1}^{n}p_{i,k_{i}}  &  =\underbrace{\sum
_{\kappa\in\operatorname*{Map}\left(  \left[  n\right]  ,\left[  m\right]
\right)  }}_{=\sum_{\kappa:\left[  n\right]  \rightarrow\left[  m\right]  }%
}\ \ \prod_{i=1}^{n}p_{i,\kappa\left(  i\right)  }\\
&  =\sum_{\kappa:\left[  n\right]  \rightarrow\left[  m\right]  }%
\ \ \prod_{i=1}^{n}p_{i,\kappa\left(  i\right)  }.
\end{align*}
Thus, (\ref{pf.lem.prodrule2.1}) becomes%
\[
\prod_{i=1}^{n}\ \ \sum_{k=1}^{m}p_{i,k}=\sum_{\left(  k_{1},k_{2}%
,\ldots,k_{n}\right)  \in\underbrace{\left[  m\right]  \times\left[  m\right]
\times\cdots\times\left[  m\right]  }_{n\text{ factors}}}\prod_{i=1}%
^{n}p_{i,k_{i}}=\sum_{\kappa:\left[  n\right]  \rightarrow\left[  m\right]
}\ \ \prod_{i=1}^{n}p_{i,\kappa\left(  i\right)  }.
\]
Lemma \ref{lem.prodrule2} is proven.
\end{proof}

Now we are ready to prove what is probably the most important property of
determinants, known as the \textit{multiplicativity of the determinant}%
\footnote{Theorem \ref{thm.det(AB)} appears, e.g., in \cite[\S 5.3, Theorem
3]{HoffmanKunze}, in \cite[proof of Theorem 5.7]{Laue-det}, in \cite[Theorem
B.17]{Strick13}, in \cite["Multiplications of determinants\textquotedblright,
Lemma]{Mate14}, in \cite[Theorem 11.4.2]{Pinkha15}, in \cite[Proposition
7.9]{GalQua18}, in \cite[\S 5]{Zeilbe}, and in \cite[Theorem 9.54]{Loehr-BC}.}:

\begin{theorem}
\label{thm.det(AB)}Let $n\in\mathbb{N}$. Let $A$ and $B$ be two $n\times
n$-matrices. Then,%
\[
\det\left(  AB\right)  =\det A\cdot\det B.
\]

\end{theorem}

\begin{proof}
[Proof of Theorem \ref{thm.det(AB)}.]Write $A$ and $B$ in the forms $A=\left(
a_{i,j}\right)  _{1\leq i\leq n,\ 1\leq j\leq n}$ and $B=\left(
b_{i,j}\right)  _{1\leq i\leq n,\ 1\leq j\leq n}$. The definition of $AB$ thus
yields $AB=\left(  \sum_{k=1}^{n}a_{i,k}b_{k,j}\right)  _{1\leq i\leq
n,\ 1\leq j\leq n}$. Therefore, (\ref{eq.det.eq.2}) (applied to $AB$ and
$\sum_{k=1}^{n}a_{i,k}b_{k,j}$ instead of $A$ and $a_{i,j}$) yields%
\begin{align}
\det\left(  AB\right)   &  =\sum_{\sigma\in S_{n}}\left(  -1\right)  ^{\sigma
}\underbrace{\prod_{i=1}^{n}\left(  \sum_{k=1}^{n}a_{i,k}b_{k,\sigma\left(
i\right)  }\right)  }_{\substack{=\sum_{\kappa:\left[  n\right]
\rightarrow\left[  n\right]  }\prod_{i=1}^{n}\left(  a_{i,\kappa\left(
i\right)  }b_{\kappa\left(  i\right)  ,\sigma\left(  i\right)  }\right)
\\\text{(by Lemma \ref{lem.prodrule2}, applied to }m=n\\\text{and }%
p_{i,k}=a_{i,k}b_{k,\sigma\left(  i\right)  }\text{)}}}\nonumber\\
&  =\sum_{\sigma\in S_{n}}\left(  -1\right)  ^{\sigma}\sum_{\kappa:\left[
n\right]  \rightarrow\left[  n\right]  }\underbrace{\prod_{i=1}^{n}\left(
a_{i,\kappa\left(  i\right)  }b_{\kappa\left(  i\right)  ,\sigma\left(
i\right)  }\right)  }_{=\left(  \prod_{i=1}^{n}a_{i,\kappa\left(  i\right)
}\right)  \left(  \prod_{i=1}^{n}b_{\kappa\left(  i\right)  ,\sigma\left(
i\right)  }\right)  }\nonumber\\
&  =\sum_{\sigma\in S_{n}}\left(  -1\right)  ^{\sigma}\sum_{\kappa:\left[
n\right]  \rightarrow\left[  n\right]  }\left(  \prod_{i=1}^{n}a_{i,\kappa
\left(  i\right)  }\right)  \left(  \prod_{i=1}^{n}b_{\kappa\left(  i\right)
,\sigma\left(  i\right)  }\right) \nonumber\\
&  =\sum_{\kappa:\left[  n\right]  \rightarrow\left[  n\right]  }\left(
\prod_{i=1}^{n}a_{i,\kappa\left(  i\right)  }\right)  \left(  \sum_{\sigma\in
S_{n}}\left(  -1\right)  ^{\sigma}\prod_{i=1}^{n}b_{\kappa\left(  i\right)
,\sigma\left(  i\right)  }\right)  . \label{pf.thm.det(AB).4}%
\end{align}

Now, for every $\kappa:\left[  n\right]  \rightarrow\left[  n\right]  $, we
let $B_{\kappa}$ be the $n\times n$-matrix $\left(  b_{\kappa\left(  i\right)
,j}\right)  _{1\leq i\leq n,\ 1\leq j\leq n}$. Then, for every $\kappa:\left[
n\right]  \rightarrow\left[  n\right]  $, the equality (\ref{eq.det.eq.2})
(applied to $B_{\kappa}$ and $b_{\kappa\left(  i\right)  ,j}$ instead of $A$
and $a_{i,j}$) yields%
\begin{equation}
\det\left(  B_{\kappa}\right)  =\sum_{\sigma\in S_{n}}\left(  -1\right)
^{\sigma}\prod_{i=1}^{n}b_{\kappa\left(  i\right)  ,\sigma\left(  i\right)  }.
\label{pf.thm.det(AB).5}%
\end{equation}
Thus, (\ref{pf.thm.det(AB).4}) becomes%
\begin{align*}
\det\left(  AB\right)   &  =\sum_{\kappa:\left[  n\right]  \rightarrow\left[
n\right]  }\left(  \prod_{i=1}^{n}a_{i,\kappa\left(  i\right)  }\right)
\underbrace{\left(  \sum_{\sigma\in S_{n}}\left(  -1\right)  ^{\sigma}%
\prod_{i=1}^{n}b_{\kappa\left(  i\right)  ,\sigma\left(  i\right)  }\right)
}_{\substack{=\det\left(  B_{\kappa}\right)  \\\text{(by
(\ref{pf.thm.det(AB).5}))}}}\\
&  =\sum_{\kappa:\left[  n\right]  \rightarrow\left[  n\right]  }\left(
\prod_{i=1}^{n}a_{i,\kappa\left(  i\right)  }\right)  \det\left(  B_{\kappa
}\right) \\
&  =\underbrace{\sum_{\substack{\kappa:\left[  n\right]  \rightarrow\left[
n\right]  ;\\\kappa\in S_{n}}}}_{\substack{=\sum_{\kappa\in S_{n}%
}\\\text{(since every }\kappa\in S_{n}\text{ automatically}\\\text{is a map
}\left[  n\right]  \rightarrow\left[  n\right]  \text{)}}}\left(  \prod
_{i=1}^{n}a_{i,\kappa\left(  i\right)  }\right)  \underbrace{\det\left(
B_{\kappa}\right)  }_{\substack{=\left(  -1\right)  ^{\kappa}\cdot\det
B\\\text{(by Lemma \ref{lem.det.sigma} \textbf{(a)})}}}\\
&  \ \ \ \ \ \ \ \ \ \ +\sum_{\substack{\kappa:\left[  n\right]
\rightarrow\left[  n\right]  ;\\\kappa\notin S_{n}}}\left(  \prod_{i=1}%
^{n}a_{i,\kappa\left(  i\right)  }\right)  \underbrace{\det\left(  B_{\kappa
}\right)  }_{\substack{=0\\\text{(by Lemma \ref{lem.det.sigma} \textbf{(b)})}%
}}\\
&  =\sum_{\kappa\in S_{n}}\left(  \prod_{i=1}^{n}a_{i,\kappa\left(  i\right)
}\right)  \left(  -1\right)  ^{\kappa}\cdot\det B+\underbrace{\sum
_{\substack{\kappa:\left[  n\right]  \rightarrow\left[  n\right]
;\\\kappa\notin S_{n}}}\left(  \prod_{i=1}^{n}a_{i,\kappa\left(  i\right)
}\right)  0}_{=0}\\
&  =\sum_{\kappa\in S_{n}}\left(  \prod_{i=1}^{n}a_{i,\kappa\left(  i\right)
}\right)  \left(  -1\right)  ^{\kappa}\cdot\det B=\sum_{\sigma\in S_{n}%
}\left(  \prod_{i=1}^{n}a_{i,\sigma\left(  i\right)  }\right)  \left(
-1\right)  ^{\sigma}\cdot\det B\\
&  \ \ \ \ \ \ \ \ \ \ \ \ \ \ \ \ \ \ \ \ \left(  \text{here, we renamed the
summation index }\kappa\text{ as }\sigma\right) \\
&  =\underbrace{\left(  \sum_{\sigma\in S_{n}}\left(  -1\right)  ^{\sigma
}\prod_{i=1}^{n}a_{i,\sigma\left(  i\right)  }\right)  }_{\substack{=\det
A\\\text{(by (\ref{eq.det.eq.2}))}}}\cdot\det B=\det A\cdot\det B.
\end{align*}
This proves Theorem \ref{thm.det(AB)}.
\end{proof}

\begin{remark}
The analogue of Theorem \ref{thm.det(AB)} with addition instead of
multiplication does not hold. If $A$ and $B$ are two $n\times n$-matrices for
some $n\in\mathbb{N}$, then $\det\left(  A+B\right)  $ does usually
\textbf{not} equal $\det A+\det B$. (However, there is a formula for
$\det\left(  A+B\right)  $ in terms of determinants of \textit{submatrices} of
$A$ and $B$; see Theorem \ref{thm.det(A+B)} below for this.)
\end{remark}

We shall now show several applications of Theorem \ref{thm.det(AB)}. First, a
simple corollary:

\begin{corollary}
\label{cor.det.product}Let $n\in\mathbb{N}$.

\textbf{(a)} If $B_{1},B_{2},\ldots,B_{k}$ are finitely many $n\times
n$-matrices, then $\det\left(  B_{1}B_{2}\cdots B_{k}\right)  =\prod_{i=1}%
^{k}\det\left(  B_{i}\right)  $.

\textbf{(b)} If $B$ is any $n\times n$-matrix, and $k\in\mathbb{N}$, then
$\det\left(  B^{k}\right)  =\left(  \det B\right)  ^{k}$.
\end{corollary}

\begin{vershort}

\begin{proof}
[Proof of Corollary \ref{cor.det.product}.]Corollary \ref{cor.det.product}
easily follows from Theorem \ref{thm.det(AB)} by induction over $k$. (The
induction base, $k=0$, relies on the fact that the product of $0$ matrices is
$I_{n}$ and has determinant $\det\left(  I_{n}\right)  =1$.) We leave the
details to the reader.
\end{proof}
\end{vershort}

\begin{verlong}

\begin{proof}
[Proof of Corollary \ref{cor.det.product}.]\textbf{(a)} Let $B_{1}%
,B_{2},\ldots,B_{k}$ be finitely many $n\times n$-matrices. We shall show that%
\begin{equation}
\det\left(  B_{1}B_{2}\cdots B_{u}\right)  =\prod_{i=1}^{u}\det\left(
B_{i}\right)  \label{pf.cor.det.product.a.1}%
\end{equation}
for every $u\in\left\{  0,1,\ldots,k\right\}  $.

We shall prove (\ref{pf.cor.det.product.a.1}) by induction over $u$:

\textit{Induction base:} We have $\det\left(  \underbrace{B_{1}B_{2}\cdots
B_{0}}_{=\left(  \text{empty product of }n\times n\text{-matrices}\right)
=I_{n}}\right)  =\det\left(  I_{n}\right)  =1$ and $\prod_{i=1}^{0}\det\left(
B_{i}\right)  =\left(  \text{empty product of elements of }\mathbb{K}\right)
=1$. Hence, $\det\left(  B_{1}B_{2}\cdots B_{0}\right)  =1=\prod_{i=1}^{0}%
\det\left(  B_{i}\right)  $. In other words, (\ref{pf.cor.det.product.a.1})
holds for $u=0$. This completes the induction base.

\textit{Induction step:} Let $U\in\left\{  0,1,\ldots,k\right\}  $ be
positive. Assume that (\ref{pf.cor.det.product.a.1}) holds for $u=U-1$. We
need to show that (\ref{pf.cor.det.product.a.1}) holds for $u=U$.

We have assumed that (\ref{pf.cor.det.product.a.1}) holds for $u=U-1$. In
other words,
\[
\det\left(  B_{1}B_{2}\cdots B_{U-1}\right)  =\prod_{i=1}^{U-1}\det\left(
B_{i}\right)  .
\]
Now,%
\begin{align*}
\det\underbrace{\left(  B_{1}B_{2}\cdots B_{U}\right)  }_{=\left(  B_{1}%
B_{2}\cdots B_{U-1}\right)  B_{U}}  &  =\det\left(  \left(  B_{1}B_{2}\cdots
B_{U-1}\right)  B_{U}\right)  =\underbrace{\det\left(  B_{1}B_{2}\cdots
B_{U-1}\right)  }_{=\prod_{i=1}^{U-1}\det\left(  B_{i}\right)  }\cdot
\det\left(  B_{U}\right) \\
&  \ \ \ \ \ \ \ \ \ \ \left(  \text{by Theorem \ref{thm.det(AB)}, applied to
}A=B_{1}B_{2}\cdots B_{U-1}\text{ and }B=B_{U}\right) \\
&  =\left(  \prod_{i=1}^{U-1}\det\left(  B_{i}\right)  \right)  \cdot
\det\left(  B_{U}\right)  =\prod_{i=1}^{U}\det\left(  B_{i}\right)  .
\end{align*}
In other words, (\ref{pf.cor.det.product.a.1}) holds for $u=U$. This completes
the induction step. Thus, (\ref{pf.cor.det.product.a.1}) is proven by induction.

Now, (\ref{pf.cor.det.product.a.1}) (applied to $u=k$) yields $\det\left(
B_{1}B_{2}\cdots B_{k}\right)  =\prod_{i=1}^{k}\det\left(  B_{i}\right)  $.
Corollary \ref{cor.det.product} \textbf{(a)} is thus proven.

\textbf{(b)} Let $B$ be any $n\times n$-matrix. Applying Corollary
\ref{cor.det.product} \textbf{(a)} to $B_{i}=B$, we obtain $\det\left(
\underbrace{BB\cdots B}_{k\text{ factors}}\right)  =\prod_{i=1}^{k}\det
B=\left(  \det B\right)  ^{k}$. Since $\underbrace{BB\cdots B}_{k\text{
factors}}=B^{k}$, this rewrites as $\det\left(  B^{k}\right)  =\left(  \det
B\right)  ^{k}$. Corollary \ref{cor.det.product} \textbf{(b)} is thus proven.
\end{proof}
\end{verlong}

\begin{example}
\label{exam.det(AB).fibo}Recall that the Fibonacci sequence is the sequence
$\left(  f_{0},f_{1},f_{2},\ldots\right)  $ of integers which is defined
recursively by $f_{0}=0$, $f_{1}=1$, and $f_{n}=f_{n-1}+f_{n-2}$ for all
$n\geq2$. We shall prove that%
\begin{equation}
f_{n+1}f_{n-1}-f_{n}^{2}=\left(  -1\right)  ^{n}\ \ \ \ \ \ \ \ \ \ \text{for
every positive integer }n. \label{eq.exam.det(AB).fibo.1}%
\end{equation}
(This is a classical fact known as
\href{https://en.wikipedia.org/wiki/Cassini_and_Catalan_identities}{the
\textit{Cassini identity}} and easy to prove by induction, but we shall prove
it differently to illustrate the use of determinants.)

Let $B$ be the $2\times2$-matrix $\left(
\begin{array}
[c]{cc}%
1 & 1\\
1 & 0
\end{array}
\right)  $ (over the ring $\mathbb{Z}$). It is easy to see that $\det B=-1$.
However, for every positive integer $n$, we have%
\begin{equation}
B^{n}=\left(
\begin{array}
[c]{cc}%
f_{n+1} & f_{n}\\
f_{n} & f_{n-1}%
\end{array}
\right)  . \label{eq.exam.det(AB).fibo.2}%
\end{equation}
Indeed, (\ref{eq.exam.det(AB).fibo.2}) can be easily proven by induction over
$n$: For $n=1$ it is clear by inspection; if it holds for $n=N$, then for
$n=N+1$ it follows from%
\begin{align*}
B^{N+1}  &  =\underbrace{B^{N}}_{\substack{=\left(
\begin{array}
[c]{cc}%
f_{N+1} & f_{N}\\
f_{N} & f_{N-1}%
\end{array}
\right)  \\\text{(by the induction hypothesis)}}}\underbrace{B}_{=\left(
\begin{array}
[c]{cc}%
1 & 1\\
1 & 0
\end{array}
\right)  }=\left(
\begin{array}
[c]{cc}%
f_{N+1} & f_{N}\\
f_{N} & f_{N-1}%
\end{array}
\right)  \left(
\begin{array}
[c]{cc}%
1 & 1\\
1 & 0
\end{array}
\right) \\
&  =\left(
\begin{array}
[c]{cc}%
f_{N+1}\cdot1+f_{N}\cdot1 & f_{N+1}\cdot1+f_{N}\cdot0\\
f_{N}\cdot1+f_{N-1}\cdot1 & f_{N}\cdot1+f_{N-1}\cdot0
\end{array}
\right) \\
&  \ \ \ \ \ \ \ \ \ \ \ \ \ \ \ \ \ \ \ \ \left(  \text{by the definition of
a product of two matrices}\right) \\
&  =\left(
\begin{array}
[c]{cc}%
f_{N+1}+f_{N} & f_{N+1}\\
f_{N}+f_{N-1} & f_{N}%
\end{array}
\right)  =\left(
\begin{array}
[c]{cc}%
f_{N+2} & f_{N+1}\\
f_{N+1} & f_{N}%
\end{array}
\right)
\end{align*}
(since $f_{N+1}+f_{N}=f_{N+2}$ and $f_{N}+f_{N-1}=f_{N+1}$).

Now, let $n$ be a positive integer. Then, (\ref{eq.exam.det(AB).fibo.2})
yields%
\[
\det\left(  B^{n}\right)  =\det\left(
\begin{array}
[c]{cc}%
f_{n+1} & f_{n}\\
f_{n} & f_{n-1}%
\end{array}
\right)  =f_{n+1}f_{n-1}-f_{n}^{2}.
\]
On the other hand, Corollary \ref{cor.det.product} \textbf{(b)} (applied to
$k=n$) yields $\det\left(  B^{n}\right)  =\left(  \underbrace{\det B}%
_{=-1}\right)  ^{n}=\left(  -1\right)  ^{n}$. Hence, $f_{n+1}f_{n-1}-f_{n}%
^{2}=\det\left(  B^{n}\right)  =\left(  -1\right)  ^{n}$. This proves
(\ref{eq.exam.det(AB).fibo.1}).
\end{example}

We can generalize (\ref{eq.exam.det(AB).fibo.1}) as follows:

\begin{exercise}
\label{exe.ps4.det.fibo}Let $a$ and $b$ be two complex numbers. Let $\left(
x_{0},x_{1},x_{2},\ldots\right)  $ be a sequence of complex numbers such that
every $n\geq2$ satisfies%
\begin{equation}
x_{n}=ax_{n-1}+bx_{n-2}. \label{eq.det.fibo.rec}%
\end{equation}
(We called such sequences \textquotedblleft$\left(  a,b\right)  $%
-recurrent\textquotedblright\ in Definition \ref{def.abrec}.) Let
$k\in\mathbb{N}$. Prove that
\begin{equation}
x_{n+1}x_{n-k-1}-x_{n}x_{n-k}=\left(  -b\right)  ^{n-k-1}\left(  x_{k+2}%
x_{0}-x_{k+1}x_{1}\right)  . \label{eq.det.fibo.claim}%
\end{equation}
for every integer $n>k$.
\end{exercise}

We notice that (\ref{eq.exam.det(AB).fibo.1}) can be obtained by applying
(\ref{eq.det.fibo.claim}) to $a=1$, $b=1$, $x_{i}=f_{i}$ and $k=0$. Thus,
(\ref{eq.det.fibo.claim}) is a generalization of (\ref{eq.exam.det(AB).fibo.1}%
). Notice that you could have easily come up with the identity
(\ref{eq.det.fibo.claim}) by trying to generalize the proof of
(\ref{eq.exam.det(AB).fibo.1}) we gave; in contrast, it is not that
straightforward to guess the general formula (\ref{eq.det.fibo.claim}) from
the classical proof of (\ref{eq.exam.det(AB).fibo.1}) by induction. So the
proof of (\ref{eq.exam.det(AB).fibo.1}) using determinants has at least the
advantage of pointing to a generalization.

\begin{example}
\label{exam.xiyj.3}Let $n\in\mathbb{N}$. Let $x_{1},x_{2},\ldots,x_{n}$ be $n$
elements of $\mathbb{K}$, and let $y_{1},y_{2},\ldots,y_{n}$ be $n$ further
elements of $\mathbb{K}$. Let $A$ be the $n\times n$-matrix $\left(
x_{i}y_{j}\right)  _{1\leq i\leq n,\ 1\leq j\leq n}$. In Example
\ref{exam.xiyj}, we have shown that $\det A=0$ if $n\geq2$. We can now prove
this in a simpler way.

Namely, let $n\geq2$. Define an $n\times n$-matrix $B$ by $B=\left(
\begin{array}
[c]{ccccc}%
x_{1} & 0 & 0 & \cdots & 0\\
x_{2} & 0 & 0 & \cdots & 0\\
x_{3} & 0 & 0 & \cdots & 0\\
\vdots & \vdots & \vdots & \ddots & \vdots\\
x_{n} & 0 & 0 & \cdots & 0
\end{array}
\right)  $. (Thus, the first column of $B$ is $\left(  x_{1},x_{2}%
,\ldots,x_{n}\right)  ^{T}$, while all other columns are filled with zeroes.)
Define an $n\times n$-matrix $C$ by $C=\left(
\begin{array}
[c]{ccccc}%
y_{1} & y_{2} & y_{3} & \cdots & y_{n}\\
0 & 0 & 0 & \cdots & 0\\
0 & 0 & 0 & \cdots & 0\\
\vdots & \vdots & \vdots & \ddots & \vdots\\
0 & 0 & 0 & \cdots & 0
\end{array}
\right)  $. (Thus, the first row of $C$ is $\left(  y_{1},y_{2},\ldots
,y_{n}\right)  $, while all other rows are filled with zeroes.)

The second row of $C$ consists of zeroes (and this second row indeed exists,
because $n\geq2$). Thus, Exercise \ref{exe.ps4.6} \textbf{(c)} (applied to $C$
instead of $A$) yields $\det C=0$. Similarly, using Exercise \ref{exe.ps4.6}
\textbf{(d)}, we can show that $\det B=0$. Now, Theorem \ref{thm.det(AB)}
(applied to $B$ and $C$ instead of $A$ and $B$) yields $\det\left(  BC\right)
=\det B\cdot\underbrace{\det C}_{=0}=0$. But what is $BC$ ?

Write $B$ in the form $B=\left(  b_{i,j}\right)  _{1\leq i\leq n,\ 1\leq j\leq
n}$, and write $C$ in the form $C=\left(  c_{i,j}\right)  _{1\leq i\leq
n,\ 1\leq j\leq n}$. Then, the definition of $BC$ yields%
\[
BC=\left(  \sum_{k=1}^{n}b_{i,k}c_{k,j}\right)  _{1\leq i\leq n,\ 1\leq j\leq
n}.
\]
Therefore, for every $\left(  i,j\right)  \in\left\{  1,2,\ldots,n\right\}
^{2}$, the $\left(  i,j\right)  $-th entry of the matrix $BC$ is%
\[
\sum_{k=1}^{n}b_{i,k}c_{k,j}=\underbrace{b_{i,1}}_{=x_{i}}\underbrace{c_{1,j}%
}_{=y_{j}}+\sum_{k=2}^{n}\underbrace{b_{i,k}}_{=0}\underbrace{c_{k,j}}%
_{=0}=x_{i}y_{j}+\underbrace{\sum_{k=2}^{n}0\cdot0}_{=0}=x_{i}y_{j}.
\]
But this is the same as the $\left(  i,j\right)  $-th entry of the matrix $A$.
Thus, every entry of $BC$ equals the corresponding entry of $A$. Hence,
$BC=A$, so that $\det\left(  BC\right)  =\det A$. Thus, $\det A=\det\left(
BC\right)  =0$, just as we wanted to show.
\end{example}

\begin{example}
\label{exam.xi+yj.3}Let $n\in\mathbb{N}$. Let $x_{1},x_{2},\ldots,x_{n}$ be
$n$ elements of $\mathbb{K}$, and let $y_{1},y_{2},\ldots,y_{n}$ be $n$
further elements of $\mathbb{K}$. Let $A$ be the $n\times n$-matrix $\left(
x_{i}+y_{j}\right)  _{1\leq i\leq n,\ 1\leq j\leq n}$. In Example
\ref{exam.xi+yj}, we have shown that $\det A=0$ if $n\geq3$.

We can now prove this in a simpler way. The argument is similar to Example
\ref{exam.xiyj.3}, and so I will be very brief:

Let $n\geq3$. Define an $n\times n$-matrix $B$ by $B=\left(
\begin{array}
[c]{ccccc}%
x_{1} & 1 & 0 & \cdots & 0\\
x_{2} & 1 & 0 & \cdots & 0\\
x_{3} & 1 & 0 & \cdots & 0\\
\vdots & \vdots & \vdots & \ddots & \vdots\\
x_{n} & 1 & 0 & \cdots & 0
\end{array}
\right)  $. (Thus, the first column of $B$ is $\left(  x_{1},x_{2}%
,\ldots,x_{n}\right)  ^{T}$, the second column is $\left(  1,1,\ldots
,1\right)  ^{T}$, while all other columns are filled with zeroes.) Define an
$n\times n$-matrix $C$ by $C=\left(
\begin{array}
[c]{ccccc}%
1 & 1 & 1 & \cdots & 1\\
y_{1} & y_{2} & y_{3} & \cdots & y_{n}\\
0 & 0 & 0 & \cdots & 0\\
\vdots & \vdots & \vdots & \ddots & \vdots\\
0 & 0 & 0 & \cdots & 0
\end{array}
\right)  $. (Thus, the first row of $C$ is $\left(  1,1,\ldots,1\right)  $,
the second row is $\left(  y_{1},y_{2},\ldots,y_{n}\right)  $, while all other
rows are filled with zeroes.) It is now easy to show that $BC=A$ (check
this!), but both $\det B$ and $\det C$ are $0$ (due to having a column or a
row filled with zeroes). Thus, again, we obtain $\det A=0$.
\end{example}

\begin{exercise}
\label{exe.ps4.pascal}Let $n\in\mathbb{N}$. Let $A$ be the $n\times n$-matrix
$\left(  \dbinom{i+j-2}{i-1}\right)  _{1\leq i\leq n,\ 1\leq j\leq n}=\left(
\begin{array}
[c]{cccc}%
\dbinom{0}{0} & \dbinom{1}{0} & \cdots & \dbinom{n-1}{0}\\
\dbinom{1}{1} & \dbinom{2}{1} & \cdots & \dbinom{n}{1}\\
\vdots & \vdots & \ddots & \vdots\\
\dbinom{n-1}{n-1} & \dbinom{n}{n-1} & \cdots & \dbinom{2n-2}{n-1}%
\end{array}
\right)  $. (This matrix is a piece of Pascal's triangle \textquotedblleft
rotated by $45^{\circ}$\textquotedblright. For example, for $n=4$, we have
$A=\left(
\begin{array}
[c]{cccc}%
1 & 1 & 1 & 1\\
1 & 2 & 3 & 4\\
1 & 3 & 6 & 10\\
1 & 4 & 10 & 20
\end{array}
\right)  $.)

Show that $\det A=1$.
\end{exercise}

The matrix $A$ in Exercise \ref{exe.ps4.pascal} is one of the so-called
\textit{Pascal matrices}; see \cite{EdelStrang} for an enlightening exposition
of some of its properties (but beware of the fact that the very first page
reveals a significant part of the solution of Exercise \ref{exe.ps4.pascal}).

\begin{remark}
There exists a more general notion of a matrix, in which the rows and the
columns are indexed not necessarily by integers from $1$ to $n$ (for some
$n\in\mathbb{N}$), but rather by arbitrary objects. For instance, this more
general notion allows us to speak of a matrix with two rows labelled
\textquotedblleft spam\textquotedblright\ and \textquotedblleft
eggs\textquotedblright, and with three columns labelled $0$, $3$ and $\infty$.
(It thus has $6$ entries, such as the $\left(  \text{\textquotedblleft
spam\textquotedblright},3\right)  $-th entry or the $\left(
\text{\textquotedblleft eggs\textquotedblright},\infty\right)  $-th entry.)
This notion of matrices is more general and more flexible than the one used
above (e.g., it allows for infinite matrices), although it has some drawbacks
(e.g., notions such as \textquotedblleft lower-triangular\textquotedblright%
\ are not defined per se, because there might be no canonical way to order the
rows and the columns; also, infinite matrices cannot always be multiplied). We
might want to define the determinant of such a matrix. Of course, this only
makes sense when the rows of the matrix are indexed by the same objects as its
columns (this essentially says that the matrix is a \textquotedblleft square
matrix\textquotedblright\ in a reasonably general sense). So, let $X$ be a
set, and $A$ be a \textquotedblleft generalized matrix\textquotedblright%
\ whose rows and columns are both indexed by the elements of $X$. We want to
define $\det A$. We assume that $X$ is finite (indeed, while $\det A$
sometimes makes sense for infinite $X$, this only happens under some rather
restrictive conditions). Then, we can define $\det A$ by
\[
\det A=\sum_{\sigma\in S_{X}}\left(  -1\right)  ^{\sigma}\prod_{i\in
X}a_{i,\sigma\left(  i\right)  },
\]
where $S_{X}$ denotes the set of all permutations of $X$. This relies on a
definition of $\left(  -1\right)  ^{\sigma}$ for every $\sigma\in S_{X}$;
fortunately, we have provided such a definition in Exercise \ref{exe.ps4.2}.
\end{remark}

We shall see more about determinants later. So far we have barely scratched
the surface. Huge collections of problems and examples on the computation of
determinants can be found in \cite{Prasolov} and \cite{Krattenthaler} (and, if
you can be bothered with 100-years-old notation and level of rigor, in Muir's
five-volume book series \cite{Muir} -- one of the most comprehensive
collections of \textquotedblleft forgotten tales\textquotedblright\ in
mathematics\footnote{In this series, Muir endeavors to summarize every paper
that had been written about determinants until the year 1920. Several of these
papers contain results that have fallen into oblivion, and not always justly
so; Muir's summaries are thus a goldmine of interesting material. However, his
notation is antiquated and his exposition is often extremely unintelligible
(e.g., complicated identities are often presented by showing an example and
hoping that the reader will correctly guess the pattern; others are stated in
verbose sentences spanning multiple lines); very few proofs are given.
\par
Three other classical British texts on determinants are Muir's and Metzler's
\cite{MuiMet60}, Turnbull's \cite{Turnbu29} and Aitken's \cite{Aitken56};
these texts (particularly the first two) contain a wealth of remarkable
results, many of which are barely remembered today. Unfortunately, their
clarity and their level of rigor leave much to be desired by modern standards.
The two-volume treatise \cite{Vodick50} and \cite{Vodick51} by Vodicka (in
Czech) might be one of the most modern texts in this classical tradition.}).

Let us finish this section with a brief remark on the geometrical use of determinants.

\begin{remark}
Let us consider the Euclidean plane $\mathbb{R}^{2}$ with its Cartesian
coordinate system and its origin $0$. If $A=\left(  x_{A},y_{A}\right)  $ and
$B=\left(  x_{B},y_{B}\right)  $ are two points on $\mathbb{R}^{2}$, then the
area of the triangle $0AB$ is $\dfrac{1}{2}\left\vert \det\left(
\begin{array}
[c]{cc}%
x_{A} & x_{B}\\
y_{A} & y_{B}%
\end{array}
\right)  \right\vert $. The absolute value here reflects the fact that
determinants can be negative, while areas must always be $\geq0$ (although
they can be $0$ when $0$, $A$ and $B$ are collinear); however, it makes
working with areas somewhat awkward. This can be circumvented by the notion of
a \textit{signed area}. (The signed area of a triangle $ABC$ is its regular
area if the triangle is \textquotedblleft directed clockwise\textquotedblright%
, and otherwise it is the negative of its area.) The signed area of the
triangle $0AB$ is $\dfrac{1}{2}\det\left(
\begin{array}
[c]{cc}%
x_{A} & x_{B}\\
y_{A} & y_{B}%
\end{array}
\right)  $. (Here, $0$ stands for the origin; i.e., \textquotedblleft the
triangle $0AB$\textquotedblright\ means the triangle with vertices at the
origin, at $A$ and at $B$.)

If $A=\left(  x_{A},y_{A}\right)  $, $B=\left(  x_{B},y_{B}\right)  $ and
$C=\left(  x_{C},y_{C}\right)  $ are three points in $\mathbb{R}^{2}$, then
the signed area of triangle $ABC$ is $\dfrac{1}{2}\det\left(
\begin{array}
[c]{ccc}%
x_{A} & x_{B} & x_{C}\\
y_{A} & y_{B} & y_{C}\\
1 & 1 & 1
\end{array}
\right)  $.

Similar formulas hold for tetrahedra: If $A=\left(  x_{A},y_{A},z_{A}\right)
$, $B=\left(  x_{B},y_{B},z_{B}\right)  $ and $C=\left(  x_{C},y_{C}%
,z_{C}\right)  $ are three points in $\mathbb{R}^{3}$, then the signed volume
of the tetrahedron $0ABC$ is $\dfrac{1}{6}\det\left(
\begin{array}
[c]{ccc}%
x_{A} & x_{B} & x_{C}\\
y_{A} & y_{B} & y_{C}\\
z_{A} & z_{B} & z_{C}%
\end{array}
\right)  $. (Again, take the absolute value for the non-signed volume.) There
is a $4\times4$ determinant formula for the signed volume of a general
tetrahedron $ABCD$.

More generally, formulas like this hold in the $n$-dimensional space
$\mathbb{R}^{n}$. Indeed, the notion of a triangle in $\mathbb{R}^{2}$ and the
notion of a tetrahedron in $\mathbb{R}^{3}$ can be generalized to a notion of
a \textit{simplex} in $\mathbb{R}^{n}$. The signed volume of such simplices
can then be defined using determinants. Then, one can extend this definition
to arbitrary polytopes (roughly speaking, convex bodies with flat bounding
faces, as opposed to \textquotedblleft round\textquotedblright\ ones like
spheres). While the general definition of the volume of a convex body uses
analysis (specifically, integrals and Lebesgue measure), such a purely
algebraic definition has its own advantages.
\end{remark}

\subsection{\label{sect.det.CB}The Cauchy-Binet formula}

This section is devoted to the Cauchy-Binet formula: a generalization of
Theorem \ref{thm.det(AB)} which is less well-known than the latter, but still
comes useful. This formula appears in the literature in various forms; we
follow \href{http://planetmath.org/cauchybinetformula}{the one on PlanetMath}
(although we use different notations).

First, we introduce a notation for \textquotedblleft picking out some rows of
a matrix and throwing away the rest\textquotedblright\ (and also the analogous
thing for columns):

\begin{definition}
\label{def.rowscols}Let $n\in\mathbb{N}$ and $m\in\mathbb{N}$. Let $A=\left(
a_{i,j}\right)  _{1\leq i\leq n,\ 1\leq j\leq m}$ be an $n\times m$-matrix.

\textbf{(a)} If $i_{1},i_{2},\ldots,i_{u}$ are some elements of $\left\{
1,2,\ldots,n\right\}  $, then we let $\operatorname*{rows}\nolimits_{i_{1}%
,i_{2},\ldots,i_{u}}A$ denote the $u\times m$-matrix $\left(  a_{i_{x}%
,j}\right)  _{1\leq x\leq u,\ 1\leq j\leq m}$. For instance, if $A=\left(
\begin{array}
[c]{ccc}%
a & a^{\prime} & a^{\prime\prime}\\
b & b^{\prime} & b^{\prime\prime}\\
c & c^{\prime} & c^{\prime\prime}\\
d & d^{\prime} & d^{\prime\prime}%
\end{array}
\right)  $, then $\operatorname*{rows}\nolimits_{3,1,4}A=\left(
\begin{array}
[c]{ccc}%
c & c^{\prime} & c^{\prime\prime}\\
a & a^{\prime} & a^{\prime\prime}\\
d & d^{\prime} & d^{\prime\prime}%
\end{array}
\right)  $. For every $p\in\left\{  1,2,\ldots,u\right\}  $, we have%
\begin{align}
&  \left(  \text{the }p\text{-th row of }\operatorname*{rows}\nolimits_{i_{1}%
,i_{2},\ldots,i_{u}}A\right) \nonumber\\
&  =\left(  a_{i_{p},1},a_{i_{p},2},\ldots,a_{i_{p},m}\right)
\ \ \ \ \ \ \ \ \ \ \left(  \text{since }\operatorname*{rows}\nolimits_{i_{1}%
,i_{2},\ldots,i_{u}}A=\left(  a_{i_{x},j}\right)  _{1\leq x\leq u,\ 1\leq
j\leq m}\right) \nonumber\\
&  =\left(  \text{the }i_{p}\text{-th row of }A\right)
\ \ \ \ \ \ \ \ \ \ \left(  \text{since }A=\left(  a_{i,j}\right)  _{1\leq
i\leq n,\ 1\leq j\leq m}\right)  . \label{eq.def.rowscols.a.row=row}%
\end{align}
Thus, $\operatorname*{rows}\nolimits_{i_{1},i_{2},\ldots,i_{u}}A$ is the
$u\times m$-matrix whose rows (from top to bottom) are the rows labelled
$i_{1},i_{2},\ldots,i_{u}$ of the matrix $A$.

\textbf{(b)} If $j_{1},j_{2},\ldots,j_{v}$ are some elements of $\left\{
1,2,\ldots,m\right\}  $, then we let $\operatorname*{cols}\nolimits_{j_{1}%
,j_{2},\ldots,j_{v}}A$ denote the $n\times v$-matrix $\left(  a_{i,j_{y}%
}\right)  _{1\leq i\leq n,\ 1\leq y\leq v}$. For instance, if $A=\left(
\begin{array}
[c]{ccc}%
a & a^{\prime} & a^{\prime\prime}\\
b & b^{\prime} & b^{\prime\prime}\\
c & c^{\prime} & c^{\prime\prime}%
\end{array}
\right)  $, then $\operatorname*{cols}\nolimits_{3,2}A=\left(
\begin{array}
[c]{cc}%
a^{\prime\prime} & a^{\prime}\\
b^{\prime\prime} & b^{\prime}\\
c^{\prime\prime} & c^{\prime}%
\end{array}
\right)  $. For every $q\in\left\{  1,2,\ldots,v\right\}  $, we have%
\begin{align}
&  \left(  \text{the }q\text{-th column of }\operatorname*{cols}%
\nolimits_{j_{1},j_{2},\ldots,j_{v}}A\right) \nonumber\\
&  =\left(
\begin{array}
[c]{c}%
a_{1,j_{q}}\\
a_{2,j_{q}}\\
\vdots\\
a_{n,j_{q}}%
\end{array}
\right)  \ \ \ \ \ \ \ \ \ \ \left(  \text{since }\operatorname*{cols}%
\nolimits_{j_{1},j_{2},\ldots,j_{v}}A=\left(  a_{i,j_{y}}\right)  _{1\leq
i\leq n,\ 1\leq y\leq v}\right) \nonumber\\
&  =\left(  \text{the }j_{q}\text{-th column of }A\right)
\ \ \ \ \ \ \ \ \ \ \left(  \text{since }A=\left(  a_{i,j}\right)  _{1\leq
i\leq n,\ 1\leq j\leq m}\right)  . \label{eq.def.rowscols.b.col=col}%
\end{align}
Thus, $\operatorname*{cols}\nolimits_{j_{1},j_{2},\ldots,j_{v}}A$ is the
$n\times v$-matrix whose columns (from left to right) are the columns labelled
$j_{1},j_{2},\ldots,j_{v}$ of the matrix $A$.
\end{definition}

Now we can state the \textit{Cauchy-Binet formula}:

\begin{theorem}
\label{thm.cauchy-binet}Let $n\in\mathbb{N}$ and $m\in\mathbb{N}$. Let $A$ be
an $n\times m$-matrix, and let $B$ be an $m\times n$-matrix. Then,%
\begin{equation}
\det\left(  AB\right)  =\sum_{1\leq g_{1}<g_{2}<\cdots<g_{n}\leq m}\det\left(
\operatorname*{cols}\nolimits_{g_{1},g_{2},\ldots,g_{n}}A\right)  \cdot
\det\left(  \operatorname*{rows}\nolimits_{g_{1},g_{2},\ldots,g_{n}}B\right)
. \label{eq.thm.cauchy-binet}%
\end{equation}

\end{theorem}

\begin{remark}
\label{rmk.cauchy-binet.sumsign}The summation sign $\sum_{1\leq g_{1}%
<g_{2}<\cdots<g_{n}\leq m}$ in (\ref{eq.thm.cauchy-binet}) is an abbreviation
for
\begin{equation}
\sum_{\substack{\left(  g_{1},g_{2},\ldots,g_{n}\right)  \in\left\{
1,2,\ldots,m\right\}  ^{n};\\g_{1}<g_{2}<\cdots<g_{n}}}.
\label{eq.rmk.cauchy-binet.real-meaning}%
\end{equation}
In particular, if $n=0$, then it signifies a summation over all $0$-tuples of
elements of $\left\{  1,2,\ldots,m\right\}  $ (because in this case, the chain
of inequalities $g_{1}<g_{2}<\cdots<g_{n}$ is a tautology); such a sum always
has exactly one addend (because there is exactly one $0$-tuple).

When both $n$ and $m$ equal $0$, then the notation $\sum_{1\leq g_{1}%
<g_{2}<\cdots<g_{n}\leq m}$ is slightly confusing: It appears to mean an empty
summation (because $1\leq m$ does not hold). But as we said, we mean this
notation to be an abbreviation for (\ref{eq.rmk.cauchy-binet.real-meaning}),
which signifies a sum with exactly one addend. But this is enough pedantry for
now; for $n>0$, the notation $\sum_{1\leq g_{1}<g_{2}<\cdots<g_{n}\leq m}$
fortunately means exactly what it seems to mean.
\end{remark}

We shall soon give a detailed proof of Theorem \ref{thm.cauchy-binet}; see
\cite[Chapter 32, Theorem]{AigZie} for a different proof\footnote{Note that
the formulation of Theorem \ref{thm.cauchy-binet} in \cite[Chapter 32,
Theorem]{AigZie} is slightly different: In our notations, it says that if $A$
is an $n\times m$-matrix and if $B$ is an $m\times n$-matrix, then%
\begin{equation}
\det\left(  AB\right)  =\sum_{\substack{\mathcal{Z}\subseteq\left\{
1,2,\ldots,m\right\}  ;\\\left\vert \mathcal{Z}\right\vert =n}}\det\left(
\operatorname*{cols}\nolimits_{\mathcal{Z}}A\right)  \cdot\det\left(
\operatorname*{rows}\nolimits_{\mathcal{Z}}B\right)  ,
\label{eq.thm.cauchy-binet.fn1.alt}%
\end{equation}
where the matrices $\operatorname*{cols}\nolimits_{\mathcal{Z}}A$ and
$\operatorname*{rows}\nolimits_{\mathcal{Z}}B$ (for $\mathcal{Z}$ being a
subset of $\left\{  1,2,\ldots,m\right\}  $) are defined as follows: Write the
subset $\mathcal{Z}$ in the form $\left\{  z_{1},z_{2},\ldots,z_{k}\right\}  $
with $z_{1}<z_{2}<\cdots<z_{k}$, and set $\operatorname*{cols}%
\nolimits_{\mathcal{Z}}A=\operatorname*{cols}\nolimits_{z_{1},z_{2}%
,\ldots,z_{k}}A$ and $\operatorname*{rows}\nolimits_{\mathcal{Z}%
}B=\operatorname*{rows}\nolimits_{z_{1},z_{2},\ldots,z_{k}}B$. (Apart from
this, \cite[Chapter 32, Theorem]{AigZie} also requires $n\leq m$; but this
requirement is useless.)
\par
The equalities (\ref{eq.thm.cauchy-binet}) and
(\ref{eq.thm.cauchy-binet.fn1.alt}) are equivalent, because the $n$-tuples
$\left(  g_{1},g_{2},\ldots,g_{n}\right)  \in\left\{  1,2,\ldots,m\right\}
^{n}$ satisfying $g_{1}<g_{2}<\cdots<g_{n}$ are in a bijection with the
subsets $\mathcal{Z}$ of $\left\{  1,2,\ldots,m\right\}  $ satisfying
$\left\vert \mathcal{Z}\right\vert =n$. (This bijection sends an $n$-tuple
$\left(  g_{1},g_{2},\ldots,g_{n}\right)  $ to the subset $\left\{
g_{1},g_{2},\ldots,g_{n}\right\}  $.)
\par
The proof of (\ref{eq.thm.cauchy-binet.fn1.alt}) given in \cite[Chapter
32]{AigZie} uses the \textit{Lindstr\"{o}m-Gessel-Viennot lemma} (which it
calls the \textquotedblleft lemma of Gessel-Viennot\textquotedblright) and is
highly worth reading.
\par
A closely related combinatorial proof of Theorem \ref{thm.cauchy-binet}
appears in \cite[\S 2]{Zeng93}. Two other proofs (one of which is more or less
the one we shall give below) are sketched in \cite[Theorem 2.3 and Problem
28.7]{Prasolov}. Yet another proof (using characteristic polynomials and block
matrices) can be found in \cite[(APP.1.2)]{LLPT95}.}.

Before we prove Theorem \ref{thm.cauchy-binet}, let us give some examples for
its use. First, here is a simple fact:

\begin{lemma}
\label{lem.increasing-sequences}Let $n\in\mathbb{N}$.

\textbf{(a)} There exists exactly one $n$-tuple $\left(  g_{1},g_{2}%
,\ldots,g_{n}\right)  \in\left\{  1,2,\ldots,n\right\}  ^{n}$ satisfying
$g_{1}<g_{2}<\cdots<g_{n}$, namely the $n$-tuple $\left(  1,2,\ldots,n\right)
$.

\textbf{(b)} Let $m\in\mathbb{N}$ be such that $m<n$. Then, there exists no
$n$-tuple $\left(  g_{1},g_{2},\ldots,g_{n}\right)  \in\left\{  1,2,\ldots
,m\right\}  ^{n}$ satisfying $g_{1}<g_{2}<\cdots<g_{n}$.
\end{lemma}

As for its intuitive meaning, Lemma \ref{lem.increasing-sequences} can be
viewed as a \textquotedblleft pigeonhole principle\textquotedblright\ for
strictly increasing sequences: Part \textbf{(b)} says (roughly speaking) that
there is no way to squeeze a strictly increasing sequence $\left(  g_{1}%
,g_{2},\ldots,g_{n}\right)  $ of $n$ numbers into the set $\left\{
1,2,\ldots,m\right\}  $ when $m<n$; part \textbf{(a)} says (again, informally)
that the only such sequence for $m=n$ is $\left(  1,2,\ldots,n\right)  $.

\begin{exercise}
\label{exe.lem.increasing-sequences} Give a formal proof of Lemma
\ref{lem.increasing-sequences}. (Do not bother doing this if you do not
particularly care about formal proofs and find Lemma
\ref{lem.increasing-sequences} obvious enough.)
\end{exercise}

\begin{example}
\label{exam.cauchy-binet.nxn}Let $n\in\mathbb{N}$. Let $A$ and $B$ be two
$n\times n$-matrices. It is easy to check that $\operatorname*{cols}%
\nolimits_{1,2,\ldots,n}A=A$ and $\operatorname*{rows}\nolimits_{1,2,\ldots
,n}B=B$. Now, Theorem \ref{thm.cauchy-binet} (applied to $m=n$) yields%
\begin{equation}
\det\left(  AB\right)  =\sum_{1\leq g_{1}<g_{2}<\cdots<g_{n}\leq n}\det\left(
\operatorname*{cols}\nolimits_{g_{1},g_{2},\ldots,g_{n}}A\right)  \cdot
\det\left(  \operatorname*{rows}\nolimits_{g_{1},g_{2},\ldots,g_{n}}B\right)
. \label{eq.exam.cauchy-binet.nxn.1}%
\end{equation}
But Lemma \ref{lem.increasing-sequences} \textbf{(a)} yields that there exists
exactly one $n$-tuple $\left(  g_{1},g_{2},\ldots,g_{n}\right)  \in\left\{
1,2,\ldots,n\right\}  ^{n}$ satisfying $g_{1}<g_{2}<\cdots<g_{n}$, namely the
$n$-tuple $\left(  1,2,\ldots,n\right)  $. Hence, the sum on the right hand
side of (\ref{eq.exam.cauchy-binet.nxn.1}) has exactly one addend: namely, the
addend for $\left(  g_{1},g_{2},\ldots,g_{n}\right)  =\left(  1,2,\ldots
,n\right)  $. Therefore, this sum simplifies as follows:%
\begin{align*}
&  \sum_{1\leq g_{1}<g_{2}<\cdots<g_{n}\leq n}\det\left(  \operatorname*{cols}%
\nolimits_{g_{1},g_{2},\ldots,g_{n}}A\right)  \cdot\det\left(
\operatorname*{rows}\nolimits_{g_{1},g_{2},\ldots,g_{n}}B\right) \\
&  =\det\left(  \underbrace{\operatorname*{cols}\nolimits_{1,2,\ldots,n}%
A}_{=A}\right)  \cdot\det\left(  \underbrace{\operatorname*{rows}%
\nolimits_{1,2,\ldots,n}B}_{=B}\right)  =\det A\cdot\det B.
\end{align*}
Hence, (\ref{eq.exam.cauchy-binet.nxn.1}) becomes%
\begin{align*}
\det\left(  AB\right)   &  =\sum_{1\leq g_{1}<g_{2}<\cdots<g_{n}\leq n}%
\det\left(  \operatorname*{cols}\nolimits_{g_{1},g_{2},\ldots,g_{n}}A\right)
\cdot\det\left(  \operatorname*{rows}\nolimits_{g_{1},g_{2},\ldots,g_{n}%
}B\right) \\
&  =\det A\cdot\det B.
\end{align*}
This, of course, is the statement of Theorem \ref{thm.det(AB)}. Hence, Theorem
\ref{thm.det(AB)} is a particular case of Theorem \ref{thm.cauchy-binet}.
\end{example}

\begin{example}
\label{exam.cauchy-binet.0}Let $n\in\mathbb{N}$ and $m\in\mathbb{N}$ be such
that $m<n$. Thus, Lemma \ref{lem.increasing-sequences} \textbf{(b)} shows that
there exists no $n$-tuple $\left(  g_{1},g_{2},\ldots,g_{n}\right)
\in\left\{  1,2,\ldots,m\right\}  ^{n}$ satisfying $g_{1}<g_{2}<\cdots<g_{n}$.

Now, let $A$ be an $n\times m$-matrix, and let $B$ be an $m\times n$-matrix.
Then, Theorem \ref{thm.cauchy-binet} yields%
\begin{align}
\det\left(  AB\right)   &  =\sum_{1\leq g_{1}<g_{2}<\cdots<g_{n}\leq m}%
\det\left(  \operatorname*{cols}\nolimits_{g_{1},g_{2},\ldots,g_{n}}A\right)
\cdot\det\left(  \operatorname*{rows}\nolimits_{g_{1},g_{2},\ldots,g_{n}%
}B\right) \nonumber\\
&  =\left(  \text{empty sum}\right) \nonumber\\
&  \ \ \ \ \ \ \ \ \ \ \left(
\begin{array}
[c]{c}%
\text{since there exists no }n\text{-tuple }\left(  g_{1},g_{2},\ldots
,g_{n}\right)  \in\left\{  1,2,\ldots,m\right\}  ^{n}\\
\text{satisfying }g_{1}<g_{2}<\cdots<g_{n}%
\end{array}
\right) \nonumber\\
&  =0. \label{eq.exam.cauchy-binet.0}%
\end{align}
This identity allows us to compute $\det A$ in Example \ref{exam.xiyj.3} in a
simpler way: Instead of defining two $n\times n$-matrices $B$ and $C$ by
$B=\left(
\begin{array}
[c]{ccccc}%
x_{1} & 0 & 0 & \cdots & 0\\
x_{2} & 0 & 0 & \cdots & 0\\
x_{3} & 0 & 0 & \cdots & 0\\
\vdots & \vdots & \vdots & \ddots & \vdots\\
x_{n} & 0 & 0 & \cdots & 0
\end{array}
\right)  $ and $C=\left(
\begin{array}
[c]{ccccc}%
y_{1} & y_{2} & y_{3} & \cdots & y_{n}\\
0 & 0 & 0 & \cdots & 0\\
0 & 0 & 0 & \cdots & 0\\
\vdots & \vdots & \vdots & \ddots & \vdots\\
0 & 0 & 0 & \cdots & 0
\end{array}
\right)  $, it suffices to define an $n\times1$-matrix $B^{\prime}$ by
$B^{\prime}=\left(  x_{1},x_{2},\ldots,x_{n}\right)  ^{T}$ and a $1\times
n$-matrix $C^{\prime}$ by $C^{\prime}=\left(  y_{1},y_{2},\ldots,y_{n}\right)
$, and argue that $A=B^{\prime}C^{\prime}$. (We leave the details to the
reader.) Similarly, Example \ref{exam.xi+yj.3} could be dealt with.
\end{example}

\begin{remark}
The equality (\ref{eq.exam.cauchy-binet.0}) can also be derived from Theorem
\ref{thm.det(AB)}. Indeed, let $n\in\mathbb{N}$ and $m\in\mathbb{N}$ be such
that $m<n$. Let $A$ be an $n\times m$-matrix, and let $B$ be an $m\times
n$-matrix. Notice that $n-m>0$ (since $m<n$). Let $A^{\prime}$ be the $n\times
n$-matrix obtained from $A$ by appending $n-m$ new columns to the right of $A$
and filling these columns with zeroes. (For example, if $n=4$ and $m=2$ and
$A=\left(
\begin{array}
[c]{cc}%
a_{1,1} & a_{1,2}\\
a_{2,1} & a_{2,2}\\
a_{3,1} & a_{3,2}\\
a_{4,1} & a_{4,2}%
\end{array}
\right)  $, then $A^{\prime}=\left(
\begin{array}
[c]{cccc}%
a_{1,1} & a_{1,2} & 0 & 0\\
a_{2,1} & a_{2,2} & 0 & 0\\
a_{3,1} & a_{3,2} & 0 & 0\\
a_{4,1} & a_{4,2} & 0 & 0
\end{array}
\right)  $.) Also, let $B^{\prime}$ be the $n\times n$-matrix obtained from
$B$ by appending $n-m$ new rows to the bottom of $B$ and filling these rows
with zeroes. (For example, if $n=4$ and $m=2$ and $B=\left(
\begin{array}
[c]{cccc}%
b_{1,1} & b_{1,2} & b_{1,3} & b_{1,4}\\
b_{2,1} & b_{2,2} & b_{2,3} & b_{2,4}%
\end{array}
\right)  $, then $B^{\prime}=\left(
\begin{array}
[c]{cccc}%
b_{1,1} & b_{1,2} & b_{1,3} & b_{1,4}\\
b_{2,1} & b_{2,2} & b_{2,3} & b_{2,4}\\
0 & 0 & 0 & 0\\
0 & 0 & 0 & 0
\end{array}
\right)  $.) Then, it is easy to check that $AB=A^{\prime}B^{\prime}$ (in
fact, just compare corresponding entries of $AB$ and $A^{\prime}B^{\prime}$).
But recall that $n-m>0$. Hence, the matrix $A^{\prime}$ has a column
consisting of zeroes (namely, its last column). Thus, Exercise \ref{exe.ps4.6}
\textbf{(d)} (applied to $A^{\prime}$ instead of $A$) shows that $\det\left(
A^{\prime}\right)  =0$. Now,
\begin{align*}
\det\left(  \underbrace{AB}_{=A^{\prime}B^{\prime}}\right)   &  =\det\left(
A^{\prime}B^{\prime}\right)  =\underbrace{\det\left(  A^{\prime}\right)
}_{=0}\cdot\det\left(  B^{\prime}\right) \\
&  \ \ \ \ \ \ \ \ \ \ \left(  \text{by Theorem \ref{thm.det(AB)}, applied to
}A^{\prime}\text{ and }B^{\prime}\text{ instead of }A\text{ and }B\right) \\
&  =0.
\end{align*}
Thus, (\ref{eq.exam.cauchy-binet.0}) is proven again.
\end{remark}

\begin{example}
\label{exam.cauchy-binet.1}Let us see what Theorem \ref{thm.cauchy-binet} says
for $n=1$. Indeed, let $m\in\mathbb{N}$; let $A=\left(  a_{1},a_{2}%
,\ldots,a_{m}\right)  $ be a $1\times m$-matrix (i.e., a row vector with $m$
entries), and let $B=\left(  b_{1},b_{2},\ldots,b_{m}\right)  ^{T}$ be an
$m\times1$-matrix (i.e., a column vector with $m$ entries). Then, $AB$ is the
$1\times1$-matrix $\left(
\begin{array}
[c]{c}%
\sum_{k=1}^{m}a_{k}b_{k}%
\end{array}
\right)  $. Thus,
\begin{equation}
\det\left(  AB\right)  =\det\left(
\begin{array}
[c]{c}%
\sum_{k=1}^{m}a_{k}b_{k}%
\end{array}
\right)  =\sum_{k=1}^{m}a_{k}b_{k}\ \ \ \ \ \ \ \ \ \ \left(  \text{by
(\ref{eq.det.small.1x1})}\right)  . \label{eq.exam.cauchy-binet.1}%
\end{equation}
What would we obtain if we tried to compute $\det\left(  AB\right)  $ using
Theorem \ref{thm.cauchy-binet}? Theorem \ref{thm.cauchy-binet} (applied to
$n=1$) yields%
\begin{align*}
\det\left(  AB\right)   &  =\underbrace{\sum_{1\leq g_{1}\leq m}}%
_{=\sum_{g_{1}=1}^{m}}\det\left(  \underbrace{\operatorname*{cols}%
\nolimits_{g_{1}}A}_{=\left(
\begin{array}
[c]{c}%
a_{g_{1}}%
\end{array}
\right)  }\right)  \cdot\det\left(  \underbrace{\operatorname*{rows}%
\nolimits_{g_{1}}B}_{=\left(
\begin{array}
[c]{c}%
b_{g_{1}}%
\end{array}
\right)  }\right) \\
&  =\sum_{g_{1}=1}^{m}\underbrace{\det\left(
\begin{array}
[c]{c}%
a_{g_{1}}%
\end{array}
\right)  }_{\substack{=a_{g_{1}}\\\text{(by (\ref{eq.det.small.1x1}))}}%
}\cdot\underbrace{\det\left(
\begin{array}
[c]{c}%
b_{g_{1}}%
\end{array}
\right)  }_{\substack{=b_{g_{1}}\\\text{(by (\ref{eq.det.small.1x1}))}}%
}=\sum_{g_{1}=1}^{m}a_{g_{1}}\cdot b_{g_{1}}.
\end{align*}
This is, of course, the same result as that of (\ref{eq.exam.cauchy-binet.1})
(with the summation index $k$ renamed as $g_{1}$). So we did not gain any
interesting insight from applying Theorem \ref{thm.cauchy-binet} to $n=1$.
\end{example}

\begin{example}
\label{exam.cauchy-binet.2}Let us try a slightly less trivial case. Indeed,
let $m\in\mathbb{N}$; let $A=\left(
\begin{array}
[c]{cccc}%
a_{1} & a_{2} & \cdots & a_{m}\\
a_{1}^{\prime} & a_{2}^{\prime} & \cdots & a_{m}^{\prime}%
\end{array}
\right)  $ be a $2\times m$-matrix, and let $B=\left(
\begin{array}
[c]{cc}%
b_{1} & b_{1}^{\prime}\\
b_{2} & b_{2}^{\prime}\\
\vdots & \vdots\\
b_{m} & b_{m}^{\prime}%
\end{array}
\right)  $ be an $m\times2$-matrix. Then, $AB$ is the $2\times2$-matrix
$\left(
\begin{array}
[c]{cc}%
\sum_{k=1}^{m}a_{k}b_{k} & \sum_{k=1}^{m}a_{k}b_{k}^{\prime}\\
\sum_{k=1}^{m}a_{k}^{\prime}b_{k} & \sum_{k=1}^{m}a_{k}^{\prime}b_{k}^{\prime}%
\end{array}
\right)  $. Hence,%
\begin{align}
\det\left(  AB\right)   &  =\det\left(
\begin{array}
[c]{cc}%
\sum_{k=1}^{m}a_{k}b_{k} & \sum_{k=1}^{m}a_{k}b_{k}^{\prime}\\
\sum_{k=1}^{m}a_{k}^{\prime}b_{k} & \sum_{k=1}^{m}a_{k}^{\prime}b_{k}^{\prime}%
\end{array}
\right) \nonumber\\
&  =\left(  \sum_{k=1}^{m}a_{k}b_{k}\right)  \left(  \sum_{k=1}^{m}%
a_{k}^{\prime}b_{k}^{\prime}\right)  -\left(  \sum_{k=1}^{m}a_{k}^{\prime
}b_{k}\right)  \left(  \sum_{k=1}^{m}a_{k}b_{k}^{\prime}\right)  .
\label{eq.exam.cauchy-binet.2.1}%
\end{align}

On the other hand, Theorem \ref{thm.cauchy-binet} (now applied to $n=2$)
yields
\begin{align*}
\det\left(  AB\right)   &  =\sum_{1\leq g_{1}<g_{2}\leq m}\det\left(
\operatorname*{cols}\nolimits_{g_{1},g_{2}}A\right)  \cdot\det\left(
\operatorname*{rows}\nolimits_{g_{1},g_{2}}B\right) \\
&  =\sum_{1\leq i<j\leq m}\det\left(  \underbrace{\operatorname*{cols}%
\nolimits_{i,j}A}_{=\left(
\begin{array}
[c]{cc}%
a_{i} & a_{j}\\
a_{i}^{\prime} & a_{j}^{\prime}%
\end{array}
\right)  }\right)  \cdot\det\left(  \underbrace{\operatorname*{rows}%
\nolimits_{i,j}B}_{=\left(
\begin{array}
[c]{cc}%
b_{i} & b_{i}^{\prime}\\
b_{j} & b_{j}^{\prime}%
\end{array}
\right)  }\right) \\
&  \ \ \ \ \ \ \ \ \ \ \left(
\begin{array}
[c]{c}%
\text{here, we renamed the summation indices }g_{1}\text{ and }g_{2}\\
\text{as }i\text{ and }j\text{, since double subscripts are annoying}%
\end{array}
\right) \\
&  =\sum_{1\leq i<j\leq m}\underbrace{\det\left(
\begin{array}
[c]{cc}%
a_{i} & a_{j}\\
a_{i}^{\prime} & a_{j}^{\prime}%
\end{array}
\right)  }_{=a_{i}a_{j}^{\prime}-a_{j}a_{i}^{\prime}}\cdot\underbrace{\det
\left(
\begin{array}
[c]{cc}%
b_{i} & b_{i}^{\prime}\\
b_{j} & b_{j}^{\prime}%
\end{array}
\right)  }_{=b_{i}b_{j}^{\prime}-b_{j}b_{i}^{\prime}}\\
&  =\sum_{1\leq i<j\leq m}\left(  a_{i}a_{j}^{\prime}-a_{j}a_{i}^{\prime
}\right)  \cdot\left(  b_{i}b_{j}^{\prime}-b_{j}b_{i}^{\prime}\right)  .
\end{align*}
Compared with (\ref{eq.exam.cauchy-binet.2.1}), this yields%
\begin{align}
&  \left(  \sum_{k=1}^{m}a_{k}b_{k}\right)  \left(  \sum_{k=1}^{m}%
a_{k}^{\prime}b_{k}^{\prime}\right)  -\left(  \sum_{k=1}^{m}a_{k}^{\prime
}b_{k}\right)  \left(  \sum_{k=1}^{m}a_{k}b_{k}^{\prime}\right) \nonumber\\
&  =\sum_{1\leq i<j\leq m}\left(  a_{i}a_{j}^{\prime}-a_{j}a_{i}^{\prime
}\right)  \cdot\left(  b_{i}b_{j}^{\prime}-b_{j}b_{i}^{\prime}\right)  .
\label{eq.exam.cauchy-binet.2.binet-cauchy}%
\end{align}
This identity is called
\href{https://en.wikipedia.org/wiki/Binet-Cauchy_identity}{the
\textit{Binet-Cauchy identity}} (I am not kidding -- look it up on the
Wikipedia). It is fairly easy to prove by direct computation; thus, using
Theorem \ref{thm.cauchy-binet} to prove it was quite an overkill. However,
(\ref{eq.exam.cauchy-binet.2.binet-cauchy}) might not be very easy to come up
with, whereas deriving it from Theorem \ref{thm.cauchy-binet} is
straightforward. (And Theorem \ref{thm.cauchy-binet} is easier to memorize
than (\ref{eq.exam.cauchy-binet.2.binet-cauchy}).)

Here is a neat application of (\ref{eq.exam.cauchy-binet.2.binet-cauchy}): If
$a_{1},a_{2},\ldots,a_{m}$ and $a_{1}^{\prime},a_{2}^{\prime},\ldots
,a_{m}^{\prime}$ are real numbers, then
(\ref{eq.exam.cauchy-binet.2.binet-cauchy}) (applied to $b_{k}=a_{k}$ and
$b_{k}^{\prime}=a_{k}^{\prime}$) yields%
\begin{align*}
&  \left(  \sum_{k=1}^{m}a_{k}a_{k}\right)  \left(  \sum_{k=1}^{m}%
a_{k}^{\prime}a_{k}^{\prime}\right)  -\left(  \sum_{k=1}^{m}a_{k}^{\prime
}a_{k}\right)  \left(  \sum_{k=1}^{m}a_{k}a_{k}^{\prime}\right) \\
&  =\sum_{1\leq i<j\leq m}\underbrace{\left(  a_{i}a_{j}^{\prime}-a_{j}%
a_{i}^{\prime}\right)  \cdot\left(  a_{i}a_{j}^{\prime}-a_{j}a_{i}^{\prime
}\right)  }_{=\left(  a_{i}a_{j}^{\prime}-a_{j}a_{i}^{\prime}\right)  ^{2}%
\geq0}\geq\sum_{1\leq i<j\leq m}0=0,
\end{align*}
so that%
\[
\left(  \sum_{k=1}^{m}a_{k}a_{k}\right)  \left(  \sum_{k=1}^{m}a_{k}^{\prime
}a_{k}^{\prime}\right)  \geq\left(  \sum_{k=1}^{m}a_{k}^{\prime}a_{k}\right)
\left(  \sum_{k=1}^{m}a_{k}a_{k}^{\prime}\right)  .
\]
In other words,%
\[
\left(  \sum_{k=1}^{m}a_{k}^{2}\right)  \left(  \sum_{k=1}^{m}\left(
a_{k}^{\prime}\right)  ^{2}\right)  \geq\left(  \sum_{k=1}^{m}a_{k}%
a_{k}^{\prime}\right)  ^{2}.
\]
This is the famous
\href{https://en.wikipedia.org/wiki/Cauchy-Schwarz_inequality}{Cauchy-Schwarz
inequality}.
\end{example}

Let us now prepare for the proof of Theorem \ref{thm.cauchy-binet}. First
comes a fact which should be fairly clear:

\begin{proposition}
\label{prop.sorting}Let $n\in\mathbb{N}$. Let $a_{1},a_{2},\ldots,a_{n}$ be
$n$ integers.

\textbf{(a)} There exists a permutation $\sigma\in S_{n}$ such that
$a_{\sigma\left(  1\right)  }\leq a_{\sigma\left(  2\right)  }\leq\cdots\leq
a_{\sigma\left(  n\right)  }$.

\textbf{(b)} If $\sigma\in S_{n}$ is such that $a_{\sigma\left(  1\right)
}\leq a_{\sigma\left(  2\right)  }\leq\cdots\leq a_{\sigma\left(  n\right)  }%
$, then, for every $i\in\left\{  1,2,\ldots,n\right\}  $, the value
$a_{\sigma\left(  i\right)  }$ depends only on $a_{1},a_{2},\ldots,a_{n}$ and
$i$ (but not on $\sigma$).

\textbf{(c)} Assume that the integers $a_{1},a_{2},\ldots,a_{n}$ are distinct.
Then, there is a \textbf{unique} permutation $\sigma\in S_{n}$ such that
$a_{\sigma\left(  1\right)  }<a_{\sigma\left(  2\right)  }<\cdots
<a_{\sigma\left(  n\right)  }$.
\end{proposition}

Let me explain why this proposition should be intuitively
obvious.\footnote{See the solution of Exercise \ref{exe.sorting.basics}
further below for a formal proof of this proposition.} Proposition
\ref{prop.sorting} \textbf{(a)} says that any list $\left(  a_{1},a_{2}%
,\ldots,a_{n}\right)  $ of $n$ integers can be sorted in weakly increasing
order by means of a permutation $\sigma\in S_{n}$. Proposition
\ref{prop.sorting} \textbf{(b)} says that the result of this sorting process
is independent of how the sorting happened (although the permutation $\sigma$
will sometimes be non-unique). Proposition \ref{prop.sorting} \textbf{(c)}
says that if the integers $a_{1},a_{2},\ldots,a_{n}$ are distinct, then the
permutation $\sigma\in S_{n}$ which sorts the list $\left(  a_{1},a_{2}%
,\ldots,a_{n}\right)  $ in increasing order is uniquely determined as well. We
required $a_{1},a_{2},\ldots,a_{n}$ to be $n$ integers for the sake of
simplicity, but we could just as well have required them to be elements of any
\textit{totally ordered set} (i.e., any set with a less-than relation
satisfying some standard axioms).

The next fact looks slightly scary, but is still rather simple:

\begin{lemma}
\label{lem.cauchy-binet.EI}For every $n\in\mathbb{N}$, let $\left[  n\right]
$ denote the set $\left\{  1,2,\ldots,n\right\}  $.

Let $n\in\mathbb{N}$ and $m\in\mathbb{N}$. We let $\mathbf{E}$ be the subset%
\[
\left\{  \left(  k_{1},k_{2},\ldots,k_{n}\right)  \in\left[  m\right]
^{n}\ \mid\ \text{the integers }k_{1},k_{2},\ldots,k_{n}\text{ are
distinct}\right\}
\]
of $\left[  m\right]  ^{n}$. We let $\mathbf{I}$ be the subset%
\[
\left\{  \left(  k_{1},k_{2},\ldots,k_{n}\right)  \in\left[  m\right]
^{n}\ \mid\ k_{1}<k_{2}<\cdots<k_{n}\right\}
\]
of $\left[  m\right]  ^{n}$. Then, the map%
\begin{align*}
\mathbf{I}\times S_{n}  &  \rightarrow\mathbf{E},\\
\left(  \left(  g_{1},g_{2},\ldots,g_{n}\right)  ,\sigma\right)   &
\mapsto\left(  g_{\sigma\left(  1\right)  },g_{\sigma\left(  2\right)
},\ldots,g_{\sigma\left(  n\right)  }\right)
\end{align*}
is well-defined and is a bijection.
\end{lemma}

The intuition for Lemma \ref{lem.cauchy-binet.EI} is that every $n$-tuple of
distinct elements of $\left\{  1,2,\ldots,m\right\}  $ can be represented
uniquely as a permuted version of a strictly increasing\footnote{An $n$-tuple
$\left(  k_{1},k_{2},\ldots,k_{n}\right)  $ is said to be \textit{strictly
increasing} if and only if $k_{1}<k_{2}<\cdots<k_{n}$.} $n$-tuple of elements
of $\left\{  1,2,\ldots,m\right\}  $, and therefore, specifying an $n$-tuple
of distinct elements of $\left\{  1,2,\ldots,m\right\}  $ is tantamount to
specifying a strictly increasing $n$-tuple of elements of $\left\{
1,2,\ldots,m\right\}  $ and a permutation $\sigma\in S_{n}$ which says how
this $n$-tuple is to be permuted.\footnote{For instance, the $4$-tuple
$\left(  4,1,6,2\right)  $ of distinct elements of $\left\{  1,2,\ldots
,7\right\}  $ can be specified by specifying the strictly increasing $4$-tuple
$\left(  1,2,4,6\right)  $ (which is its sorted version) and the permutation
$\pi\in S_{4}$ which sends $1,2,3,4$ to $3,1,4,2$, respectively (that is,
$\pi=\left(  3,1,4,2\right)  $ in one-line notation). In the terminology of
Lemma \ref{lem.cauchy-binet.EI}, the map
\begin{align*}
\mathbf{I}\times S_{n}  &  \rightarrow\mathbf{E},\\
\left(  \left(  g_{1},g_{2},\ldots,g_{n}\right)  ,\sigma\right)   &
\mapsto\left(  g_{\sigma\left(  1\right)  },g_{\sigma\left(  2\right)
},\ldots,g_{\sigma\left(  n\right)  }\right)
\end{align*}
sends $\left(  \left(  1,2,4,6\right)  ,\pi\right)  $ to $\left(
4,1,6,2\right)  $.} This is not a formal proof, but this should explain why
Lemma \ref{lem.cauchy-binet.EI} is usually applied throughout mathematics
without even mentioning it as a statement. If desired, a formal proof of Lemma
\ref{lem.cauchy-binet.EI} can be obtained using Proposition \ref{prop.sorting}%
.\footnote{Again, see the solution of Exercise \ref{exe.sorting.basics}
further below for such a proof.}

\begin{exercise}
\label{exe.sorting.basics}Prove Proposition \ref{prop.sorting} and Lemma
\ref{lem.cauchy-binet.EI}. (Ignore this exercise if you find these two facts
sufficiently obvious and are uninterested in the details of their proofs.)
\end{exercise}

Before we return to Theorem \ref{thm.cauchy-binet}, let me make a digression
about sorting:

\begin{exercise}
\label{exe.sorting.nmu}Let $n\in\mathbb{N}$ and $m\in\mathbb{N}$ be such that
$n\geq m$. Let $a_{1},a_{2},\ldots,a_{n}$ be $n$ integers. Let $b_{1}%
,b_{2},\ldots,b_{m}$ be $m$ integers. Assume that%
\begin{equation}
a_{i}\leq b_{i}\ \ \ \ \ \ \ \ \ \ \text{for every }i\in\left\{
1,2,\ldots,m\right\}  . \label{eq.exe.sorting.nmu.ass}%
\end{equation}

Let $\sigma\in S_{n}$ be such that $a_{\sigma\left(  1\right)  }\leq
a_{\sigma\left(  2\right)  }\leq\cdots\leq a_{\sigma\left(  n\right)  }$. Let
$\tau\in S_{m}$ be such that $b_{\tau\left(  1\right)  }\leq b_{\tau\left(
2\right)  }\leq\cdots\leq b_{\tau\left(  m\right)  }$. Then,%
\[
a_{\sigma\left(  i\right)  }\leq b_{\tau\left(  i\right)  }%
\ \ \ \ \ \ \ \ \ \ \text{for every }i\in\left\{  1,2,\ldots,m\right\}  .
\]

\end{exercise}

\begin{remark}
\label{rmk.sorting.nmu.interpret}Loosely speaking, Exercise
\ref{exe.sorting.nmu} says the following: If two lists $\left(  a_{1}%
,a_{2},\ldots,a_{n}\right)  $ and $\left(  b_{1},b_{2},\ldots,b_{m}\right)  $
of integers have the property that each entry of the first list is $\leq$ to
the corresponding entry of the second list (as long as the latter is
well-defined), then this property still holds after both lists are sorted in
increasing order, provided that we have $n\geq m$ (that is, the first list is
at least as long as the second list).

A consequence of Exercise \ref{exe.sorting.nmu} is the following curious fact,
known as the \textquotedblleft non-messing-up phenomenon\textquotedblright%
\ (\cite[Theorem 1]{Tenner-NMU} and \cite[Example 1]{GalKar71}): If we start
with a matrix filled with integers, then sort the entries of each row of the
matrix in increasing order, and then sort the entries of each column of the
resulting matrix in increasing order, then the final matrix still has sorted
rows (i.e., the entries of each row are still sorted). That is, the sorting of
the columns did not \textquotedblleft mess up\textquotedblright\ the
sortedness of the rows. For example, if we start with the matrix $\left(
\begin{array}
[c]{cccc}%
1 & 3 & 2 & 5\\
2 & 1 & 4 & 2\\
3 & 1 & 6 & 0
\end{array}
\right)  $, then sorting the entries of each row gives us the matrix $\left(
\begin{array}
[c]{cccc}%
1 & 2 & 3 & 5\\
1 & 2 & 2 & 4\\
0 & 1 & 3 & 6
\end{array}
\right)  $, and then sorting the entries of each column results in the matrix
$\left(
\begin{array}
[c]{cccc}%
0 & 1 & 2 & 4\\
1 & 2 & 3 & 5\\
1 & 2 & 3 & 6
\end{array}
\right)  $. The rows of this matrix are still sorted, as the \textquotedblleft
non-messing-up phenomenon\textquotedblright\ predicts. To prove this
phenomenon in general, it suffices to show that any entry in the resulting
matrix is $\leq$ to the entry directly below it (assuming that the latter
exists); but this follows easily from Exercise \ref{exe.sorting.nmu}.
\end{remark}

We are now ready to prove Theorem \ref{thm.cauchy-binet}.

\begin{vershort}

\begin{proof}
[Proof of Theorem \ref{thm.cauchy-binet}.]We shall use the notations of Lemma
\ref{lem.cauchy-binet.EI}.

Write the $n\times m$-matrix $A$ as $A=\left(  a_{i,j}\right)  _{1\leq i\leq
n,\ 1\leq j\leq m}$. Write the $m\times n$-matrix $B$ as $B=\left(
b_{i,j}\right)  _{1\leq i\leq m,\ 1\leq j\leq n}$. The definition of $AB$ thus
yields $AB=\left(  \sum_{k=1}^{m}a_{i,k}b_{k,j}\right)  _{1\leq i\leq
n,\ 1\leq j\leq n}$. Therefore, (\ref{eq.det.eq.2}) (applied to $AB$ and
$\sum_{k=1}^{m}a_{i,k}b_{k,j}$ instead of $A$ and $a_{i,j}$) yields%
\begin{equation}
\det\left(  AB\right)  =\sum_{\sigma\in S_{n}}\left(  -1\right)  ^{\sigma
}\prod_{i=1}^{n}\left(  \sum_{k=1}^{m}a_{i,k}b_{k,\sigma\left(  i\right)
}\right)  . \label{pf.thm.cauchy-binet.short.1}%
\end{equation}
But for every $\sigma\in S_{n}$, we have%
\begin{align*}
\prod_{i=1}^{n}\left(  \sum_{k=1}^{m}a_{i,k}b_{k,\sigma\left(  i\right)
}\right)   &  =\sum_{\left(  k_{1},k_{2},\ldots,k_{n}\right)  \in\left[
m\right]  ^{n}}\underbrace{\prod_{i=1}^{n}\left(  a_{i,k_{i}}b_{k_{i}%
,\sigma\left(  i\right)  }\right)  }_{=\left(  \prod_{i=1}^{n}a_{i,k_{i}%
}\right)  \left(  \prod_{i=1}^{n}b_{k_{i},\sigma\left(  i\right)  }\right)
}\\
&  \ \ \ \ \ \ \ \ \ \ \left(  \text{by Lemma \ref{lem.prodrule}, applied to
}m_{i}=n\text{ and }p_{i,k}=a_{i,k}b_{k,\sigma\left(  i\right)  }\right) \\
&  =\sum_{\left(  k_{1},k_{2},\ldots,k_{n}\right)  \in\left[  m\right]  ^{n}%
}\left(  \prod_{i=1}^{n}a_{i,k_{i}}\right)  \left(  \prod_{i=1}^{n}%
b_{k_{i},\sigma\left(  i\right)  }\right)  .
\end{align*}
Hence, (\ref{pf.thm.cauchy-binet.short.1}) rewrites as%
\begin{align}
\det\left(  AB\right)   &  =\sum_{\sigma\in S_{n}}\left(  -1\right)  ^{\sigma
}\sum_{\left(  k_{1},k_{2},\ldots,k_{n}\right)  \in\left[  m\right]  ^{n}%
}\left(  \prod_{i=1}^{n}a_{i,k_{i}}\right)  \left(  \prod_{i=1}^{n}%
b_{k_{i},\sigma\left(  i\right)  }\right) \nonumber\\
&  =\underbrace{\sum_{\sigma\in S_{n}}\sum_{\left(  k_{1},k_{2},\ldots
,k_{n}\right)  \in\left[  m\right]  ^{n}}}_{=\sum_{\left(  k_{1},k_{2}%
,\ldots,k_{n}\right)  \in\left[  m\right]  ^{n}}\sum_{\sigma\in S_{n}}}\left(
-1\right)  ^{\sigma}\left(  \prod_{i=1}^{n}a_{i,k_{i}}\right)  \left(
\prod_{i=1}^{n}b_{k_{i},\sigma\left(  i\right)  }\right) \nonumber\\
&  =\sum_{\left(  k_{1},k_{2},\ldots,k_{n}\right)  \in\left[  m\right]  ^{n}%
}\sum_{\sigma\in S_{n}}\left(  -1\right)  ^{\sigma}\left(  \prod_{i=1}%
^{n}a_{i,k_{i}}\right)  \left(  \prod_{i=1}^{n}b_{k_{i},\sigma\left(
i\right)  }\right) \nonumber\\
&  =\sum_{\left(  k_{1},k_{2},\ldots,k_{n}\right)  \in\left[  m\right]  ^{n}%
}\left(  \prod_{i=1}^{n}a_{i,k_{i}}\right)  \left(  \sum_{\sigma\in S_{n}%
}\left(  -1\right)  ^{\sigma}\prod_{i=1}^{n}b_{k_{i},\sigma\left(  i\right)
}\right)  . \label{pf.thm.cauchy-binet.short.3}%
\end{align}
But every $\left(  k_{1},k_{2},\ldots,k_{n}\right)  \in\left[  m\right]  ^{n}$
satisfies%
\begin{equation}
\sum_{\sigma\in S_{n}}\left(  -1\right)  ^{\sigma}\prod_{i=1}^{n}%
b_{k_{i},\sigma\left(  i\right)  }=\det\left(  \operatorname*{rows}%
\nolimits_{k_{1},k_{2},\ldots,k_{n}}B\right)
\label{pf.thm.cauchy-binet.short.detrows}%
\end{equation}
\footnote{\textit{Proof.} Let $\left(  k_{1},k_{2},\ldots,k_{n}\right)
\in\left[  m\right]  ^{n}$. Recall that $B=\left(  b_{i,j}\right)  _{1\leq
i\leq m,\ 1\leq j\leq n}$. Hence, the definition of $\operatorname*{rows}%
\nolimits_{k_{1},k_{2},\ldots,k_{n}}B$ gives us%
\[
\operatorname*{rows}\nolimits_{k_{1},k_{2},\ldots,k_{n}}B=\left(  b_{k_{x}%
,j}\right)  _{1\leq x\leq n,\ 1\leq j\leq n}=\left(  b_{k_{i},j}\right)
_{1\leq i\leq n,\ 1\leq j\leq n}%
\]
(here, we renamed the index $x$ as $i$). Hence, (\ref{eq.det.eq.2}) (applied
to $\operatorname*{rows}\nolimits_{k_{1},k_{2},\ldots,k_{n}}B$ and
$b_{k_{i},j}$ instead of $A$ and $a_{i,j}$) yields%
\[
\det\left(  \operatorname*{rows}\nolimits_{k_{1},k_{2},\ldots,k_{n}}B\right)
=\sum_{\sigma\in S_{n}}\left(  -1\right)  ^{\sigma}\prod_{i=1}^{n}%
b_{k_{i},\sigma\left(  i\right)  },
\]
qed.}. Hence, (\ref{pf.thm.cauchy-binet.short.3}) becomes%
\begin{align}
\det\left(  AB\right)   &  =\sum_{\left(  k_{1},k_{2},\ldots,k_{n}\right)
\in\left[  m\right]  ^{n}}\left(  \prod_{i=1}^{n}a_{i,k_{i}}\right)
\underbrace{\left(  \sum_{\sigma\in S_{n}}\left(  -1\right)  ^{\sigma}%
\prod_{i=1}^{n}b_{k_{i},\sigma\left(  i\right)  }\right)  }_{=\det\left(
\operatorname*{rows}\nolimits_{k_{1},k_{2},\ldots,k_{n}}B\right)  }\nonumber\\
&  =\sum_{\left(  k_{1},k_{2},\ldots,k_{n}\right)  \in\left[  m\right]  ^{n}%
}\left(  \prod_{i=1}^{n}a_{i,k_{i}}\right)  \det\left(  \operatorname*{rows}%
\nolimits_{k_{1},k_{2},\ldots,k_{n}}B\right)  .
\label{pf.thm.cauchy-binet.short.4}%
\end{align}

But for every $\left(  k_{1},k_{2},\ldots,k_{n}\right)  \in\left[  m\right]
^{n}$ satisfying $\left(  k_{1},k_{2},\ldots,k_{n}\right)  \notin\mathbf{E}$,
we have%
\begin{equation}
\det\left(  \operatorname*{rows}\nolimits_{k_{1},k_{2},\ldots,k_{n}}B\right)
=0 \label{pf.thm.cauchy-binet.short.6a}%
\end{equation}
\footnote{\textit{Proof of (\ref{pf.thm.cauchy-binet.short.6a}):} Let $\left(
k_{1},k_{2},\ldots,k_{n}\right)  \in\left[  m\right]  ^{n}$ be such that
$\left(  k_{1},k_{2},\ldots,k_{n}\right)  \notin\mathbf{E}$. Then, the
integers $k_{1},k_{2},\ldots,k_{n}$ are not distinct (because $\mathbf{E}$ is
the set of all $n$-tuples in $\left[  m\right]  ^{n}$ whose entries are
distinct). Thus, there exist two distinct elements $p$ and $q$ of $\left[
n\right]  $ such that $k_{p}=k_{q}$. Consider these $p$ and $q$. But
$\operatorname*{rows}\nolimits_{k_{1},k_{2},\ldots,k_{n}}B$ is the $n\times
n$-matrix whose rows (from top to bottom) are the rows labelled $k_{1}%
,k_{2},\ldots,k_{n}$ of the matrix $B$. Since $k_{p} = k_{q}$, this shows that
the $p$-th row and the $q$-th row of the matrix $\operatorname*{rows}%
\nolimits_{k_{1},k_{2},\ldots,k_{n}}B$ are equal. Hence, the matrix
$\operatorname*{rows}\nolimits_{k_{1},k_{2},\ldots,k_{n}}B$ has two equal rows
(since $p$ and $q$ are distinct). Therefore, Exercise \ref{exe.ps4.6}
\textbf{(e)} (applied to $\operatorname*{rows}\nolimits_{k_{1},k_{2}%
,\ldots,k_{n}}B$ instead of $A$) yields $\det\left(  \operatorname*{rows}%
\nolimits_{k_{1},k_{2},\ldots,k_{n}}B\right)  =0$, qed.}. Therefore, in the
sum on the right hand side of (\ref{pf.thm.cauchy-binet.short.4}), all the
addends corresponding to $\left(  k_{1},k_{2},\ldots,k_{n}\right)  \in\left[
m\right]  ^{n}$ satisfying $\left(  k_{1},k_{2},\ldots,k_{n}\right)
\notin\mathbf{E}$ evaluate to $0$. We can therefore remove all these addends
from the sum. The remaining addends are those corresponding to $\left(
k_{1},k_{2},\ldots,k_{n}\right)  \in\mathbf{E}$. Therefore,
(\ref{pf.thm.cauchy-binet.short.4}) becomes
\begin{equation}
\det\left(  AB\right)  =\sum_{\left(  k_{1},k_{2},\ldots,k_{n}\right)
\in\mathbf{E}}\left(  \prod_{i=1}^{n}a_{i,k_{i}}\right)  \det\left(
\operatorname*{rows}\nolimits_{k_{1},k_{2},\ldots,k_{n}}B\right)  .
\label{pf.thm.cauchy-binet.short.7}%
\end{equation}

On the other hand, Lemma \ref{lem.cauchy-binet.EI} yields that the map%
\begin{align*}
\mathbf{I}\times S_{n}  &  \rightarrow\mathbf{E},\\
\left(  \left(  g_{1},g_{2},\ldots,g_{n}\right)  ,\sigma\right)   &
\mapsto\left(  g_{\sigma\left(  1\right)  },g_{\sigma\left(  2\right)
},\ldots,g_{\sigma\left(  n\right)  }\right)
\end{align*}
is well-defined and is a bijection. Hence, we can substitute $\left(
g_{\sigma\left(  1\right)  },g_{\sigma\left(  2\right)  },\ldots
,g_{\sigma\left(  n\right)  }\right)  $ for $\left(  k_{1},k_{2},\ldots
,k_{n}\right)  $ in the sum on the right hand side of
(\ref{pf.thm.cauchy-binet.short.7}). We thus obtain%
\begin{align*}
&  \sum_{\left(  k_{1},k_{2},\ldots,k_{n}\right)  \in\mathbf{E}}\left(
\prod_{i=1}^{n}a_{i,k_{i}}\right)  \det\left(  \operatorname*{rows}%
\nolimits_{k_{1},k_{2},\ldots,k_{n}}B\right) \\
&  =\sum_{\left(  \left(  g_{1},g_{2},\ldots,g_{n}\right)  ,\sigma\right)
\in\mathbf{I}\times S_{n}}\left(  \prod_{i=1}^{n}a_{i,g_{\sigma\left(
i\right)  }}\right)  \det\left(  \operatorname*{rows}\nolimits_{g_{\sigma
\left(  1\right)  },g_{\sigma\left(  2\right)  },\ldots,g_{\sigma\left(
n\right)  }}B\right)  .
\end{align*}
Thus, (\ref{pf.thm.cauchy-binet.short.7}) becomes%
\begin{align}
\det\left(  AB\right)   &  =\sum_{\left(  k_{1},k_{2},\ldots,k_{n}\right)
\in\mathbf{E}}\left(  \prod_{i=1}^{n}a_{i,k_{i}}\right)  \det\left(
\operatorname*{rows}\nolimits_{k_{1},k_{2},\ldots,k_{n}}B\right) \nonumber\\
&  =\sum_{\left(  \left(  g_{1},g_{2},\ldots,g_{n}\right)  ,\sigma\right)
\in\mathbf{I}\times S_{n}}\left(  \prod_{i=1}^{n}a_{i,g_{\sigma\left(
i\right)  }}\right)  \det\left(  \operatorname*{rows}\nolimits_{g_{\sigma
\left(  1\right)  },g_{\sigma\left(  2\right)  },\ldots,g_{\sigma\left(
n\right)  }}B\right)  . \label{pf.thm.cauchy-binet.short.9}%
\end{align}

But every $\left(  k_{1},k_{2},\ldots,k_{n}\right)  \in\left[  m\right]  ^{n}$
and every $\sigma\in S_{n}$ satisfy%
\begin{equation}
\det\left(  \operatorname*{rows}\nolimits_{k_{\sigma\left(  1\right)
},k_{\sigma\left(  2\right)  },\ldots,k_{\sigma\left(  n\right)  }}B\right)
=\left(  -1\right)  ^{\sigma}\cdot\det\left(  \operatorname*{rows}%
\nolimits_{k_{1},k_{2},\ldots,k_{n}}B\right)
\label{pf.thm.cauchy-binet.short.6b}%
\end{equation}
\footnote{\textit{Proof of (\ref{pf.thm.cauchy-binet.short.6b}):} Let $\left(
k_{1},k_{2},\ldots,k_{n}\right)  \in\left[  m\right]  ^{n}$ and $\sigma\in
S_{n}$. We have $\operatorname*{rows}\nolimits_{k_{1},k_{2},\ldots,k_{n}%
}B=\left(  b_{k_{i},j}\right)  _{1\leq i\leq n,\ 1\leq j\leq n}$ (as we have
seen in one of the previous footnotes) and $\operatorname*{rows}%
\nolimits_{k_{\sigma\left(  1\right)  },k_{\sigma\left(  2\right)  }%
,\ldots,k_{\sigma\left(  n\right)  }}B=\left(  b_{k_{\sigma\left(  i\right)
},j}\right)  _{1\leq i\leq n,\ 1\leq j\leq n}$ (for similar reasons). Hence,
we can apply Lemma \ref{lem.det.sigma} \textbf{(a)} to $\sigma$,
$\operatorname*{rows}\nolimits_{k_{1},k_{2},\ldots,k_{n}}B$, $b_{k_{i},j}$ and
$\operatorname*{rows}\nolimits_{k_{\sigma\left(  1\right)  },k_{\sigma\left(
2\right)  },\ldots,k_{\sigma\left(  n\right)  }}B$ instead of $\kappa$, $B$,
$b_{i,j}$ and $B_{\kappa}$. As a consequence, we obtain
\[
\det\left(  \operatorname*{rows}\nolimits_{k_{\sigma\left(  1\right)
},k_{\sigma\left(  2\right)  },\ldots,k_{\sigma\left(  n\right)  }}B\right)
=\left(  -1\right)  ^{\sigma}\cdot\det\left(  \operatorname*{rows}%
\nolimits_{k_{1},k_{2},\ldots,k_{n}}B\right)  .
\]
This proves (\ref{pf.thm.cauchy-binet.short.6b}).}. Hence,
(\ref{pf.thm.cauchy-binet.short.9}) becomes%
\begin{align}
&  \det\left(  AB\right) \nonumber\\
&  =\underbrace{\sum_{\left(  \left(  g_{1},g_{2},\ldots,g_{n}\right)
,\sigma\right)  \in\mathbf{I}\times S_{n}}}_{=\sum_{\left(  g_{1},g_{2}%
,\ldots,g_{n}\right)  \in\mathbf{I}}\sum_{\sigma\in S_{n}}}\left(  \prod
_{i=1}^{n}a_{i,g_{\sigma\left(  i\right)  }}\right)  \underbrace{\det\left(
\operatorname*{rows}\nolimits_{g_{\sigma\left(  1\right)  },g_{\sigma\left(
2\right)  },\ldots,g_{\sigma\left(  n\right)  }}B\right)  }%
_{\substack{=\left(  -1\right)  ^{\sigma}\cdot\det\left(  \operatorname*{rows}%
\nolimits_{g_{1},g_{2},\ldots,g_{n}}B\right)  \\\text{(by
(\ref{pf.thm.cauchy-binet.short.6b}), applied to }k_{i}=g_{i}\text{)}%
}}\nonumber\\
&  =\sum_{\left(  g_{1},g_{2},\ldots,g_{n}\right)  \in\mathbf{I}}\sum
_{\sigma\in S_{n}}\left(  \prod_{i=1}^{n}a_{i,g_{\sigma\left(  i\right)  }%
}\right)  \left(  -1\right)  ^{\sigma}\cdot\det\left(  \operatorname*{rows}%
\nolimits_{g_{1},g_{2},\ldots,g_{n}}B\right) \nonumber\\
&  =\sum_{\left(  g_{1},g_{2},\ldots,g_{n}\right)  \in\mathbf{I}}\left(
\sum_{\sigma\in S_{n}}\left(  \prod_{i=1}^{n}a_{i,g_{\sigma\left(  i\right)
}}\right)  \left(  -1\right)  ^{\sigma}\right)  \cdot\det\left(
\operatorname*{rows}\nolimits_{g_{1},g_{2},\ldots,g_{n}}B\right)  .
\label{pf.thm.cauchy-binet.short.10}%
\end{align}
But every $\left(  g_{1},g_{2},\ldots,g_{n}\right)  \in\mathbf{I}$ satisfies
$\sum_{\sigma\in S_{n}}\left(  \prod_{i=1}^{n}a_{i,g_{\sigma\left(  i\right)
}}\right)  \left(  -1\right)  ^{\sigma}=\det\left(  \operatorname*{cols}%
\nolimits_{g_{1},g_{2},\ldots,g_{n}}A\right)  $%
\ \ \ \ \footnote{\textit{Proof.} Let $\left(  g_{1},g_{2},\ldots
,g_{n}\right)  \in\mathbf{I}$. We have $A=\left(  a_{i,j}\right)  _{1\leq
i\leq n,\ 1\leq j\leq m}$. Thus, the definition of $\operatorname*{cols}%
\nolimits_{g_{1},g_{2},\ldots,g_{n}}A$ yields%
\[
\operatorname*{cols}\nolimits_{g_{1},g_{2},\ldots,g_{n}}A=\left(  a_{i,g_{y}%
}\right)  _{1\leq i\leq n,\ 1\leq y\leq n}=\left(  a_{i,g_{j}}\right)  _{1\leq
i\leq n,\ 1\leq j\leq n}%
\]
(here, we renamed the index $y$ as $j$). Hence, (\ref{eq.det.eq.2}) (applied
to $\operatorname*{cols}\nolimits_{g_{1},g_{2},\ldots,g_{n}}A$ and
$a_{i,g_{j}}$ instead of $A$ and $a_{i,j}$) yields%
\[
\det\left(  \operatorname*{cols}\nolimits_{g_{1},g_{2},\ldots,g_{n}}A\right)
=\sum_{\sigma\in S_{n}}\left(  -1\right)  ^{\sigma}\prod_{i=1}^{n}%
a_{i,g_{\sigma\left(  i\right)  }}=\sum_{\sigma\in S_{n}}\left(  \prod
_{i=1}^{n}a_{i,g_{\sigma\left(  i\right)  }}\right)  \left(  -1\right)
^{\sigma},
\]
qed.}. Hence, (\ref{pf.thm.cauchy-binet.short.10}) becomes%
\begin{align}
&  \det\left(  AB\right) \nonumber\\
&  =\sum_{\left(  g_{1},g_{2},\ldots,g_{n}\right)  \in\mathbf{I}%
}\underbrace{\left(  \sum_{\sigma\in S_{n}}\left(  \prod_{i=1}^{n}%
a_{i,g_{\sigma\left(  i\right)  }}\right)  \left(  -1\right)  ^{\sigma
}\right)  }_{=\det\left(  \operatorname*{cols}\nolimits_{g_{1},g_{2}%
,\ldots,g_{n}}A\right)  }\cdot\det\left(  \operatorname*{rows}\nolimits_{g_{1}%
,g_{2},\ldots,g_{n}}B\right) \nonumber\\
&  =\sum_{\left(  g_{1},g_{2},\ldots,g_{n}\right)  \in\mathbf{I}}\det\left(
\operatorname*{cols}\nolimits_{g_{1},g_{2},\ldots,g_{n}}A\right)  \cdot
\det\left(  \operatorname*{rows}\nolimits_{g_{1},g_{2},\ldots,g_{n}}B\right)
. \label{pf.thm.cauchy-binet.short.15}%
\end{align}
Finally, we recall that $\mathbf{I}$ was defined as the set
\[
\left\{  \left(  k_{1},k_{2},\ldots,k_{n}\right)  \in\left[  m\right]
^{n}\ \mid\ k_{1}<k_{2}<\cdots<k_{n}\right\}  .
\]
Thus, summing over all $\left(  g_{1},g_{2},\ldots,g_{n}\right)  \in
\mathbf{I}$ means the same as summing over all $\left(  g_{1},g_{2}%
,\ldots,g_{n}\right)  \in\left[  m\right]  ^{n}$ satisfying $g_{1}%
<g_{2}<\cdots<g_{n}$. In other words,%
\[
\sum_{\left(  g_{1},g_{2},\ldots,g_{n}\right)  \in\mathbf{I}}=\sum
_{\substack{\left(  g_{1},g_{2},\ldots,g_{n}\right)  \in\left[  m\right]
^{n};\\g_{1}<g_{2}<\cdots<g_{n}}}=\sum_{1\leq g_{1}<g_{2}<\cdots<g_{n}\leq m}%
\]
(an equality between summation signs -- hopefully its meaning is obvious).
Hence, (\ref{pf.thm.cauchy-binet.short.15}) becomes%
\[
\det\left(  AB\right)  =\sum_{1\leq g_{1}<g_{2}<\cdots<g_{n}\leq m}\det\left(
\operatorname*{cols}\nolimits_{g_{1},g_{2},\ldots,g_{n}}A\right)  \cdot
\det\left(  \operatorname*{rows}\nolimits_{g_{1},g_{2},\ldots,g_{n}}B\right)
.
\]
This proves Theorem \ref{thm.cauchy-binet}.
\end{proof}
\end{vershort}

\begin{verlong}

\begin{proof}
[Proof of Theorem \ref{thm.cauchy-binet}.]We shall use the notations of Lemma
\ref{lem.cauchy-binet.EI}. Recall that%
\[
\mathbf{E}=\left\{  \left(  k_{1},k_{2},\ldots,k_{n}\right)  \in\left[
m\right]  ^{n}\ \mid\ \text{the integers }k_{1},k_{2},\ldots,k_{n}\text{ are
distinct}\right\}
\]
and%
\begin{align*}
\mathbf{I}  &  =\left\{  \left(  k_{1},k_{2},\ldots,k_{n}\right)  \in\left[
m\right]  ^{n}\ \mid\ k_{1}<k_{2}<\cdots<k_{n}\right\} \\
&  =\left\{  \left(  k_{1},k_{2},\ldots,k_{n}\right)  \in\left\{
1,2,\ldots,m\right\}  ^{n}\ \mid\ k_{1}<k_{2}<\cdots<k_{n}\right\} \\
&  \ \ \ \ \ \ \ \ \ \ \left(  \text{since }\left[  m\right]  =\left\{
1,2,\ldots,m\right\}  \right) \\
&  =\left\{  \left(  g_{1},g_{2},\ldots,g_{n}\right)  \in\left\{
1,2,\ldots,m\right\}  ^{n}\ \mid\ g_{1}<g_{2}<\cdots<g_{n}\right\}
\end{align*}
(here, we renamed the index $\left(  k_{1},k_{2},\ldots,k_{n}\right)  $ as
$\left(  g_{1},g_{2},\ldots,g_{n}\right)  $).

Lemma \ref{lem.cauchy-binet.EI} says that the map%
\begin{align*}
\mathbf{I}\times S_{n}  &  \rightarrow\mathbf{E},\\
\left(  \left(  g_{1},g_{2},\ldots,g_{n}\right)  ,\sigma\right)   &
\mapsto\left(  g_{\sigma\left(  1\right)  },g_{\sigma\left(  2\right)
},\ldots,g_{\sigma\left(  n\right)  }\right)
\end{align*}
is well-defined and is a bijection.

Write the $n\times m$-matrix $A$ as $A=\left(  a_{i,j}\right)  _{1\leq i\leq
n,\ 1\leq j\leq m}$. Write the $m\times n$-matrix $B$ as $B=\left(
b_{i,j}\right)  _{1\leq i\leq m,\ 1\leq j\leq n}$. The definition of $AB$ thus
yields $AB=\left(  \sum_{k=1}^{m}a_{i,k}b_{k,j}\right)  _{1\leq i\leq
n,\ 1\leq j\leq n}$. Therefore, (\ref{eq.det.eq.2}) (applied to $AB$ and
$\sum_{k=1}^{m}a_{i,k}b_{k,j}$ instead of $A$ and $a_{i,j}$) yields%
\begin{equation}
\det\left(  AB\right)  =\sum_{\sigma\in S_{n}}\left(  -1\right)  ^{\sigma
}\prod_{i=1}^{n}\left(  \sum_{k=1}^{m}a_{i,k}b_{k,\sigma\left(  i\right)
}\right)  . \label{pf.thm.cauchy-binet.1}%
\end{equation}
But for every $\sigma\in S_{n}$, we have%
\begin{align*}
\prod_{i=1}^{n}\left(  \sum_{k=1}^{m}a_{i,k}b_{k,\sigma\left(  i\right)
}\right)   &  =\sum_{\left(  k_{1},k_{2},\ldots,k_{n}\right)  \in
\underbrace{\left[  m\right]  \times\left[  m\right]  \times\cdots
\times\left[  m\right]  }_{n\text{ factors}}}\prod_{i=1}^{n}\left(
a_{i,k_{i}}b_{k_{i},\sigma\left(  i\right)  }\right) \\
&  \ \ \ \ \ \ \ \ \ \ \left(  \text{by Lemma \ref{lem.prodrule}, applied to
}m_{i}=n\text{ and }p_{i,k}=a_{i,k}b_{k,\sigma\left(  i\right)  }\right) \\
&  =\sum_{\left(  k_{1},k_{2},\ldots,k_{n}\right)  \in\left[  m\right]  ^{n}%
}\underbrace{\prod_{i=1}^{n}\left(  a_{i,k_{i}}b_{k_{i},\sigma\left(
i\right)  }\right)  }_{=\left(  \prod_{i=1}^{n}a_{i,k_{i}}\right)  \left(
\prod_{i=1}^{n}b_{k_{i},\sigma\left(  i\right)  }\right)  }\\
&  \ \ \ \ \ \ \ \ \ \ \left(  \text{since }\underbrace{\left[  m\right]
\times\left[  m\right]  \times\cdots\times\left[  m\right]  }_{n\text{
factors}}=\left[  m\right]  ^{n}\right) \\
&  =\sum_{\left(  k_{1},k_{2},\ldots,k_{n}\right)  \in\left[  m\right]  ^{n}%
}\left(  \prod_{i=1}^{n}a_{i,k_{i}}\right)  \left(  \prod_{i=1}^{n}%
b_{k_{i},\sigma\left(  i\right)  }\right)  .
\end{align*}
Hence, (\ref{pf.thm.cauchy-binet.1}) becomes%
\begin{align}
\det\left(  AB\right)   &  =\sum_{\sigma\in S_{n}}\left(  -1\right)  ^{\sigma
}\underbrace{\prod_{i=1}^{n}\left(  \sum_{k=1}^{m}a_{i,k}b_{k,\sigma\left(
i\right)  }\right)  }_{=\sum_{\left(  k_{1},k_{2},\ldots,k_{n}\right)
\in\left[  m\right]  ^{n}}\left(  \prod_{i=1}^{n}a_{i,k_{i}}\right)  \left(
\prod_{i=1}^{n}b_{k_{i},\sigma\left(  i\right)  }\right)  }\nonumber\\
&  =\sum_{\sigma\in S_{n}}\left(  -1\right)  ^{\sigma}\sum_{\left(
k_{1},k_{2},\ldots,k_{n}\right)  \in\left[  m\right]  ^{n}}\left(  \prod
_{i=1}^{n}a_{i,k_{i}}\right)  \left(  \prod_{i=1}^{n}b_{k_{i},\sigma\left(
i\right)  }\right) \nonumber\\
&  =\underbrace{\sum_{\sigma\in S_{n}}\sum_{\left(  k_{1},k_{2},\ldots
,k_{n}\right)  \in\left[  m\right]  ^{n}}}_{=\sum_{\left(  k_{1},k_{2}%
,\ldots,k_{n}\right)  \in\left[  m\right]  ^{n}}\sum_{\sigma\in S_{n}}}\left(
-1\right)  ^{\sigma}\left(  \prod_{i=1}^{n}a_{i,k_{i}}\right)  \left(
\prod_{i=1}^{n}b_{k_{i},\sigma\left(  i\right)  }\right) \nonumber\\
&  =\sum_{\left(  k_{1},k_{2},\ldots,k_{n}\right)  \in\left[  m\right]  ^{n}%
}\sum_{\sigma\in S_{n}}\left(  -1\right)  ^{\sigma}\left(  \prod_{i=1}%
^{n}a_{i,k_{i}}\right)  \left(  \prod_{i=1}^{n}b_{k_{i},\sigma\left(
i\right)  }\right) \nonumber\\
&  =\sum_{\left(  k_{1},k_{2},\ldots,k_{n}\right)  \in\left[  m\right]  ^{n}%
}\left(  \prod_{i=1}^{n}a_{i,k_{i}}\right)  \left(  \sum_{\sigma\in S_{n}%
}\left(  -1\right)  ^{\sigma}\prod_{i=1}^{n}b_{k_{i},\sigma\left(  i\right)
}\right)  . \label{pf.thm.cauchy-binet.3}%
\end{align}
But every $\left(  k_{1},k_{2},\ldots,k_{n}\right)  \in\left[  m\right]  ^{n}$
satisfies $\sum_{\sigma\in S_{n}}\left(  -1\right)  ^{\sigma}\prod_{i=1}%
^{n}b_{k_{i},\sigma\left(  i\right)  }=\det\left(  \operatorname*{rows}%
\nolimits_{k_{1},k_{2},\ldots,k_{n}}B\right)  $%
\ \ \ \ \footnote{\textit{Proof.} Let $\left(  k_{1},k_{2},\ldots
,k_{n}\right)  \in\left[  m\right]  ^{n}$. The definition of
$\operatorname*{rows}\nolimits_{k_{1},k_{2},\ldots,k_{n}}B$ yields
$\operatorname*{rows}\nolimits_{k_{1},k_{2},\ldots,k_{n}}B=\left(  b_{k_{x}%
,j}\right)  _{1\leq x\leq n,\ 1\leq j\leq n}$ (since $B=\left(  b_{i,j}%
\right)  _{1\leq i\leq m,\ 1\leq j\leq n}$). Thus,%
\[
\operatorname*{rows}\nolimits_{k_{1},k_{2},\ldots,k_{n}}B=\left(  b_{k_{x}%
,j}\right)  _{1\leq x\leq n,\ 1\leq j\leq n}=\left(  b_{k_{i},j}\right)
_{1\leq i\leq n,\ 1\leq j\leq n}%
\]
(here, we renamed the index $x$ as $i$). Hence, (\ref{eq.det.eq.2}) (applied
to $\operatorname*{rows}\nolimits_{k_{1},k_{2},\ldots,k_{n}}B$ and
$b_{k_{i},j}$ instead of $A$ and $a_{i,j}$) yields%
\[
\det\left(  \operatorname*{rows}\nolimits_{k_{1},k_{2},\ldots,k_{n}}B\right)
=\sum_{\sigma\in S_{n}}\left(  -1\right)  ^{\sigma}\prod_{i=1}^{n}%
b_{k_{i},\sigma\left(  i\right)  },
\]
qed.}. Hence, (\ref{pf.thm.cauchy-binet.3}) becomes%
\begin{align}
\det\left(  AB\right)   &  =\sum_{\left(  k_{1},k_{2},\ldots,k_{n}\right)
\in\left[  m\right]  ^{n}}\left(  \prod_{i=1}^{n}a_{i,k_{i}}\right)
\underbrace{\left(  \sum_{\sigma\in S_{n}}\left(  -1\right)  ^{\sigma}%
\prod_{i=1}^{n}b_{k_{i},\sigma\left(  i\right)  }\right)  }_{=\det\left(
\operatorname*{rows}\nolimits_{k_{1},k_{2},\ldots,k_{n}}B\right)  }\nonumber\\
&  =\sum_{\left(  k_{1},k_{2},\ldots,k_{n}\right)  \in\left[  m\right]  ^{n}%
}\left(  \prod_{i=1}^{n}a_{i,k_{i}}\right)  \det\left(  \operatorname*{rows}%
\nolimits_{k_{1},k_{2},\ldots,k_{n}}B\right)  . \label{pf.thm.cauchy-binet.4}%
\end{align}

Next, let us make two observations:

\begin{itemize}
\item Every $\left(  k_{1},k_{2},\ldots,k_{n}\right)  \in\left[  m\right]
^{n}$ satisfying $\left(  k_{1},k_{2},\ldots,k_{n}\right)  \notin\mathbf{E}$
satisfies%
\begin{equation}
\det\left(  \operatorname*{rows}\nolimits_{k_{1},k_{2},\ldots,k_{n}}B\right)
=0 \label{pf.thm.cauchy-binet.6a}%
\end{equation}
\footnote{\textit{Proof of (\ref{pf.thm.cauchy-binet.6a}):} Let $\left(
g_{1},g_{2},\ldots,g_{n}\right)  \in\left[  m\right]  ^{n}$ be such that
$\left(  g_{1},g_{2},\ldots,g_{n}\right)  \notin\mathbf{E}$.
\par
Let us first show that there exist two distinct elements $p$ and $q$ of
$\left[  n\right]  $ such that $g_{p}=g_{q}$. Indeed, we assume the contrary.
Thus, there do not exist two distinct elements $p$ and $q$ of $\left[
n\right]  $ such that $g_{p}=g_{q}$. In other words, every two distinct
elements $p$ and $q$ of $\left[  n\right]  $ satisfy $g_{p}\neq g_{q}$. In
other words, the integers $g_{1},g_{2},\ldots,g_{n}$ are distinct. Hence,
$\left(  g_{1},g_{2},\ldots,g_{n}\right)  $ is an element $\left(  k_{1}%
,k_{2},\ldots,k_{n}\right)  \in\left[  m\right]  ^{n}$ such that the integers
$k_{1},k_{2},\ldots,k_{n}$ are distinct (since $\left(  g_{1},g_{2}%
,\ldots,g_{n}\right)  \in\left[  m\right]  ^{n}$ and since the integers
$g_{1},g_{2},\ldots,g_{n}$ are distinct). In other words,%
\[
\left(  g_{1},g_{2},\ldots,g_{n}\right)  \in\left\{  \left(  k_{1}%
,k_{2},\ldots,k_{n}\right)  \in\left[  m\right]  ^{n}\ \mid\ \text{the
integers }k_{1},k_{2},\ldots,k_{n}\text{ are distinct}\right\}  =\mathbf{E}.
\]
This contradicts $\left(  g_{1},g_{2},\ldots,g_{n}\right)  \notin\mathbf{E}$.
This contradiction proves that our assumption was wrong. Hence, we have proven
that there exist two distinct elements $p$ and $q$ of $\left[  n\right]  $
such that $g_{p}=g_{q}$. Consider these $p$ and $q$. Now,
(\ref{eq.def.rowscols.a.row=row}) (applied to $m$, $n$, $n$, $g_{r}$, $B$ and
$b_{i,j}$ instead of $n$, $m$, $u$, $i_{r}$, $A$ and $a_{i,j}$) yields%
\begin{align}
&  \left(  \text{the }p\text{-th row of }\operatorname*{rows}\nolimits_{g_{1}%
,g_{2},\ldots,g_{n}}B\right) \nonumber\\
&  =\left(  \text{the }\underbrace{g_{p}}_{=g_{q}}\text{-th row of }B\right)
=\left(  \text{the }g_{q}\text{-th row of }B\right)  .
\label{pf.thm.cauchy-binet.6a.pf.1}%
\end{align}
On the other hand, (\ref{eq.def.rowscols.a.row=row}) (applied to $m$, $n$,
$n$, $g_{r}$, $B$, $b_{i,j}$ and $q$ instead of $n$, $m$, $u$, $i_{r}$, $A$,
$a_{i,j}$ and $p$) yields%
\begin{equation}
\left(  \text{the }q\text{-th row of }\operatorname*{rows}\nolimits_{g_{1}%
,g_{2},\ldots,g_{n}}B\right)  =\left(  \text{the }g_{q}\text{-th row of
}B\right)  .\nonumber
\end{equation}
Compared with (\ref{pf.thm.cauchy-binet.6a.pf.1}), this yields
\[
\left(  \text{the }p\text{-th row of }\operatorname*{rows}\nolimits_{g_{1}%
,g_{2},\ldots,g_{n}}B\right)  =\left(  \text{the }q\text{-th row of
}\operatorname*{rows}\nolimits_{g_{1},g_{2},\ldots,g_{n}}B\right)  .
\]
Thus, the matrix $\operatorname*{rows}\nolimits_{g_{1},g_{2},\ldots,g_{n}}B$
has two equal rows (namely, the $p$-th row and the $q$-th row), because $p$
and $q$ are distinct. Therefore, Exercise \ref{exe.ps4.6} \textbf{(e)}
(applied to $\operatorname*{rows}\nolimits_{g_{1},g_{2},\ldots,g_{n}}B$
instead of $A$) yields $\det\left(  \operatorname*{rows}\nolimits_{g_{1}%
,g_{2},\ldots,g_{n}}B\right)  =0$.
\par
Let us now forget that we fixed $\left(  g_{1},g_{2},\ldots,g_{n}\right)  $.
We thus have shown that every $\left(  g_{1},g_{2},\ldots,g_{n}\right)
\in\left[  m\right]  ^{n}$ satisfying $\left(  g_{1},g_{2},\ldots
,g_{n}\right)  \notin\mathbf{E}$ satisfies $\det\left(  \operatorname*{rows}%
\nolimits_{g_{1},g_{2},\ldots,g_{n}}B\right)  =0$. Renaming $\left(
g_{1},g_{2},\ldots,g_{n}\right)  $ as $\left(  k_{1},k_{2},\ldots
,k_{n}\right)  $ in this result, we obtain the following: Every $\left(
k_{1},k_{2},\ldots,k_{n}\right)  \in\left[  m\right]  ^{n}$ satisfying
$\left(  k_{1},k_{2},\ldots,k_{n}\right)  \notin\mathbf{E}$ satisfies
$\det\left(  \operatorname*{rows}\nolimits_{k_{1},k_{2},\ldots,k_{n}}B\right)
=0$. This proves (\ref{pf.thm.cauchy-binet.6a}).}.

\item Every $\left(  k_{1},k_{2},\ldots,k_{n}\right)  \in\left[  m\right]
^{n}$ and every $\sigma\in S_{n}$ satisfy%
\begin{equation}
\det\left(  \operatorname*{rows}\nolimits_{k_{\sigma\left(  1\right)
},k_{\sigma\left(  2\right)  },\ldots,k_{\sigma\left(  n\right)  }}B\right)
=\left(  -1\right)  ^{\sigma}\cdot\det\left(  \operatorname*{rows}%
\nolimits_{k_{1},k_{2},\ldots,k_{n}}B\right)  \label{pf.thm.cauchy-binet.6b}%
\end{equation}
\footnote{\textit{Proof of (\ref{pf.thm.cauchy-binet.6b}):} Let $\left(
k_{1},k_{2},\ldots,k_{n}\right)  \in\left[  m\right]  ^{n}$ and $\sigma\in
S_{n}$. The definition of $\operatorname*{rows}\nolimits_{k_{1},k_{2}%
,\ldots,k_{n}}B$ yields $\operatorname*{rows}\nolimits_{k_{1},k_{2}%
,\ldots,k_{n}}B=\left(  b_{k_{x},j}\right)  _{1\leq x\leq n,\ 1\leq j\leq n}$
(since $B=\left(  b_{i,j}\right)  _{1\leq i\leq m,\ 1\leq j\leq n}$). Thus,%
\[
\operatorname*{rows}\nolimits_{k_{1},k_{2},\ldots,k_{n}}B=\left(  b_{k_{x}%
,j}\right)  _{1\leq x\leq n,\ 1\leq j\leq n}=\left(  b_{k_{i},j}\right)
_{1\leq i\leq n,\ 1\leq j\leq n}%
\]
(here, we renamed the index $x$ as $i$). On the other hand, the definition of
$\operatorname*{rows}\nolimits_{k_{\sigma\left(  1\right)  },k_{\sigma\left(
2\right)  },\ldots,k_{\sigma\left(  n\right)  }}B$ yields
$\operatorname*{rows}\nolimits_{k_{\sigma\left(  1\right)  },k_{\sigma\left(
2\right)  },\ldots,k_{\sigma\left(  n\right)  }}B=\left(  b_{k_{\sigma\left(
x\right)  },j}\right)  _{1\leq x\leq n,\ 1\leq j\leq n}$ (since $B=\left(
b_{i,j}\right)  _{1\leq i\leq m,\ 1\leq j\leq n}$). Thus,%
\[
\operatorname*{rows}\nolimits_{k_{\sigma\left(  1\right)  },k_{\sigma\left(
2\right)  },\ldots,k_{\sigma\left(  n\right)  }}B=\left(  b_{k_{\sigma\left(
x\right)  },j}\right)  _{1\leq x\leq n,\ 1\leq j\leq n}=\left(  b_{k_{\sigma
\left(  i\right)  },j}\right)  _{1\leq i\leq n,\ 1\leq j\leq n}%
\]
(here, we renamed the index $x$ as $i$). So we know that $\operatorname*{rows}%
\nolimits_{k_{1},k_{2},\ldots,k_{n}}B=\left(  b_{k_{i},j}\right)  _{1\leq
i\leq n,\ 1\leq j\leq n}$ and $\operatorname*{rows}\nolimits_{k_{\sigma\left(
1\right)  },k_{\sigma\left(  2\right)  },\ldots,k_{\sigma\left(  n\right)  }%
}B=\left(  b_{k_{\sigma\left(  i\right)  },j}\right)  _{1\leq i\leq n,\ 1\leq
j\leq n}$. Hence, we can apply Lemma \ref{lem.det.sigma} \textbf{(a)} to
$\sigma$, $\operatorname*{rows}\nolimits_{k_{1},k_{2},\ldots,k_{n}}B$,
$b_{k_{i},j}$ and $\operatorname*{rows}\nolimits_{k_{\sigma\left(  1\right)
},k_{\sigma\left(  2\right)  },\ldots,k_{\sigma\left(  n\right)  }}B$ instead
of $\kappa$, $B$, $b_{i,j}$ and $B_{\kappa}$. As a consequence, we obtain
\[
\det\left(  \operatorname*{rows}\nolimits_{k_{\sigma\left(  1\right)
},k_{\sigma\left(  2\right)  },\ldots,k_{\sigma\left(  n\right)  }}B\right)
=\left(  -1\right)  ^{\sigma}\cdot\det\left(  \operatorname*{rows}%
\nolimits_{k_{1},k_{2},\ldots,k_{n}}B\right)  .
\]
This proves (\ref{pf.thm.cauchy-binet.6b}).}
\end{itemize}

Now, (\ref{pf.thm.cauchy-binet.4}) becomes%
\begin{align*}
\det\left(  AB\right)   &  =\sum_{\left(  k_{1},k_{2},\ldots,k_{n}\right)
\in\left[  m\right]  ^{n}}\left(  \prod_{i=1}^{n}a_{i,k_{i}}\right)
\det\left(  \operatorname*{rows}\nolimits_{k_{1},k_{2},\ldots,k_{n}}B\right)
\\
&  =\underbrace{\sum_{\substack{\left(  k_{1},k_{2},\ldots,k_{n}\right)
\in\left[  m\right]  ^{n};\\\left(  k_{1},k_{2},\ldots,k_{n}\right)
\in\mathbf{E}}}}_{\substack{=\sum_{\left(  k_{1},k_{2},\ldots,k_{n}\right)
\in\mathbf{E}}\\\text{(since }\mathbf{E}\subseteq\left[  m\right]
^{n}\text{)}}}\left(  \prod_{i=1}^{n}a_{i,k_{i}}\right)  \det\left(
\operatorname*{rows}\nolimits_{k_{1},k_{2},\ldots,k_{n}}B\right) \\
&  \ \ \ \ \ \ \ \ \ \ +\sum_{\substack{\left(  k_{1},k_{2},\ldots
,k_{n}\right)  \in\left[  m\right]  ^{n};\\\left(  k_{1},k_{2},\ldots
,k_{n}\right)  \notin\mathbf{E}}}\left(  \prod_{i=1}^{n}a_{i,k_{i}}\right)
\underbrace{\det\left(  \operatorname*{rows}\nolimits_{k_{1},k_{2}%
,\ldots,k_{n}}B\right)  }_{\substack{=0\\\text{(by
(\ref{pf.thm.cauchy-binet.6a}) (since }\left(  k_{1},k_{2},\ldots
,k_{n}\right)  \notin\mathbf{E}\text{))}}}\\
&  =\sum_{\left(  k_{1},k_{2},\ldots,k_{n}\right)  \in\mathbf{E}}\left(
\prod_{i=1}^{n}a_{i,k_{i}}\right)  \det\left(  \operatorname*{rows}%
\nolimits_{k_{1},k_{2},\ldots,k_{n}}B\right) \\
&  \ \ \ \ \ \ \ \ \ \ +\underbrace{\sum_{\substack{\left(  k_{1},k_{2}%
,\ldots,k_{n}\right)  \in\left[  m\right]  ^{n};\\\left(  k_{1},k_{2}%
,\ldots,k_{n}\right)  \notin\mathbf{E}}}\left(  \prod_{i=1}^{n}a_{i,k_{i}%
}\right)  0}_{=0}\\
&  =\sum_{\left(  k_{1},k_{2},\ldots,k_{n}\right)  \in\mathbf{E}}\left(
\prod_{i=1}^{n}a_{i,k_{i}}\right)  \det\left(  \operatorname*{rows}%
\nolimits_{k_{1},k_{2},\ldots,k_{n}}B\right) \\
&  =\sum_{\left(  \left(  g_{1},g_{2},\ldots,g_{n}\right)  ,\sigma\right)
\in\mathbf{I}\times S_{n}}\left(  \prod_{i=1}^{n}a_{i,g_{\sigma\left(
i\right)  }}\right)  \det\left(  \operatorname*{rows}\nolimits_{g_{\sigma
\left(  1\right)  },g_{\sigma\left(  2\right)  },\ldots,g_{\sigma\left(
n\right)  }}B\right)
\end{align*}
(here, we have substituted $\left(  g_{\sigma\left(  1\right)  }%
,g_{\sigma\left(  2\right)  },\ldots,g_{\sigma\left(  n\right)  }\right)  $
for $\left(  k_{1},k_{2},\ldots,k_{n}\right)  $ in the sum, because the map%
\begin{align*}
\mathbf{I}\times S_{n}  &  \rightarrow\mathbf{E},\\
\left(  \left(  g_{1},g_{2},\ldots,g_{n}\right)  ,\sigma\right)   &
\mapsto\left(  g_{\sigma\left(  1\right)  },g_{\sigma\left(  2\right)
},\ldots,g_{\sigma\left(  n\right)  }\right)
\end{align*}
is a bijection). Thus,%
\begin{align}
&  \det\left(  AB\right) \nonumber\\
&  =\underbrace{\sum_{\left(  \left(  g_{1},g_{2},\ldots,g_{n}\right)
,\sigma\right)  \in\mathbf{I}\times S_{n}}}_{=\sum_{\left(  g_{1},g_{2}%
,\ldots,g_{n}\right)  \in\mathbf{I}}\sum_{\sigma\in S_{n}}}\left(  \prod
_{i=1}^{n}a_{i,g_{\sigma\left(  i\right)  }}\right)  \underbrace{\det\left(
\operatorname*{rows}\nolimits_{g_{\sigma\left(  1\right)  },g_{\sigma\left(
2\right)  },\ldots,g_{\sigma\left(  n\right)  }}B\right)  }%
_{\substack{=\left(  -1\right)  ^{\sigma}\cdot\det\left(  \operatorname*{rows}%
\nolimits_{g_{1},g_{2},\ldots,g_{n}}B\right)  \\\text{(by
(\ref{pf.thm.cauchy-binet.6b}), applied to }k_{i}=g_{i}\text{)}}}\nonumber\\
&  =\sum_{\left(  g_{1},g_{2},\ldots,g_{n}\right)  \in\mathbf{I}}\sum
_{\sigma\in S_{n}}\left(  \prod_{i=1}^{n}a_{i,g_{\sigma\left(  i\right)  }%
}\right)  \left(  -1\right)  ^{\sigma}\cdot\det\left(  \operatorname*{rows}%
\nolimits_{g_{1},g_{2},\ldots,g_{n}}B\right) \nonumber\\
&  =\sum_{\left(  g_{1},g_{2},\ldots,g_{n}\right)  \in\mathbf{I}}\left(
\sum_{\sigma\in S_{n}}\left(  \prod_{i=1}^{n}a_{i,g_{\sigma\left(  i\right)
}}\right)  \left(  -1\right)  ^{\sigma}\right)  \cdot\det\left(
\operatorname*{rows}\nolimits_{g_{1},g_{2},\ldots,g_{n}}B\right)  .
\label{pf.thm.cauchy-binet.10}%
\end{align}
But every $\left(  g_{1},g_{2},\ldots,g_{n}\right)  \in\mathbf{I}$ satisfies
$\sum_{\sigma\in S_{n}}\left(  \prod_{i=1}^{n}a_{i,g_{\sigma\left(  i\right)
}}\right)  \left(  -1\right)  ^{\sigma}=\det\left(  \operatorname*{cols}%
\nolimits_{g_{1},g_{2},\ldots,g_{n}}A\right)  $%
\ \ \ \ \footnote{\textit{Proof.} Let $\left(  g_{1},g_{2},\ldots
,g_{n}\right)  \in\mathbf{I}$. Then,%
\[
\left(  g_{1},g_{2},\ldots,g_{n}\right)  \in\mathbf{I}=\left\{  \left(
k_{1},k_{2},\ldots,k_{n}\right)  \in\left[  m\right]  ^{n}\ \mid\ k_{1}%
<k_{2}<\cdots<k_{n}\right\}  .
\]
In other words, $\left(  g_{1},g_{2},\ldots,g_{n}\right)  $ is an element
$\left(  k_{1},k_{2},\ldots,k_{n}\right)  \in\left[  m\right]  ^{n}$
satisfying $k_{1}<k_{2}<\cdots<k_{n}$. In other words, $\left(  g_{1}%
,g_{2},\ldots,g_{n}\right)  $ is an element of $\left[  m\right]  ^{n}$ and
satisfies $g_{1}<g_{2}<\cdots<g_{n}$.
\par
Now, the definition of $\operatorname*{cols}\nolimits_{g_{1},g_{2}%
,\ldots,g_{n}}A$ yields $\operatorname*{cols}\nolimits_{g_{1},g_{2}%
,\ldots,g_{n}}A=\left(  a_{i,g_{y}}\right)  _{1\leq i\leq n,\ 1\leq y\leq n}$
(since $A=\left(  a_{i,j}\right)  _{1\leq i\leq n,\ 1\leq j\leq m}$). Hence,%
\[
\operatorname*{cols}\nolimits_{g_{1},g_{2},\ldots,g_{n}}A=\left(  a_{i,g_{y}%
}\right)  _{1\leq i\leq n,\ 1\leq y\leq n}=\left(  a_{i,g_{j}}\right)  _{1\leq
i\leq n,\ 1\leq j\leq n}%
\]
(here, we renamed the index $y$ as $j$). Hence, (\ref{eq.det.eq.2}) (applied
to $\operatorname*{cols}\nolimits_{g_{1},g_{2},\ldots,g_{n}}A$ and
$a_{i,g_{j}}$ instead of $A$ and $a_{i,j}$) yields%
\[
\det\left(  \operatorname*{cols}\nolimits_{g_{1},g_{2},\ldots,g_{n}}A\right)
=\sum_{\sigma\in S_{n}}\left(  -1\right)  ^{\sigma}\prod_{i=1}^{n}%
a_{i,g_{\sigma\left(  i\right)  }}=\sum_{\sigma\in S_{n}}\left(  \prod
_{i=1}^{n}a_{i,g_{\sigma\left(  i\right)  }}\right)  \left(  -1\right)
^{\sigma},
\]
qed.}. Hence, (\ref{pf.thm.cauchy-binet.10}) becomes%
\begin{align*}
&  \det\left(  AB\right) \\
&  =\sum_{\left(  g_{1},g_{2},\ldots,g_{n}\right)  \in\mathbf{I}%
}\underbrace{\left(  \sum_{\sigma\in S_{n}}\left(  \prod_{i=1}^{n}%
a_{i,g_{\sigma\left(  i\right)  }}\right)  \left(  -1\right)  ^{\sigma
}\right)  }_{=\det\left(  \operatorname*{cols}\nolimits_{g_{1},g_{2}%
,\ldots,g_{n}}A\right)  }\cdot\det\left(  \operatorname*{rows}\nolimits_{g_{1}%
,g_{2},\ldots,g_{n}}B\right) \\
&  =\underbrace{\sum_{\left(  g_{1},g_{2},\ldots,g_{n}\right)  \in\mathbf{I}}%
}_{\substack{=\sum_{\substack{\left(  g_{1},g_{2},\ldots,g_{n}\right)
\in\left\{  1,2,\ldots,m\right\}  ^{n};\\g_{1}<g_{2}<\cdots<g_{n}%
}}\\\text{(since }\mathbf{I}=\left\{  \left(  g_{1},g_{2},\ldots,g_{n}\right)
\in\left\{  1,2,\ldots,m\right\}  ^{n}\ \mid\ g_{1}<g_{2}<\cdots
<g_{n}\right\}  \text{)}}}\det\left(  \operatorname*{cols}\nolimits_{g_{1}%
,g_{2},\ldots,g_{n}}A\right)  \cdot\det\left(  \operatorname*{rows}%
\nolimits_{g_{1},g_{2},\ldots,g_{n}}B\right) \\
&  =\underbrace{\sum_{\substack{\left(  g_{1},g_{2},\ldots,g_{n}\right)
\in\left\{  1,2,\ldots,m\right\}  ^{n};\\g_{1}<g_{2}<\cdots<g_{n}}%
}}_{\substack{=\sum_{1\leq g_{1}<g_{2}<\cdots<g_{n}\leq m}\\\text{(since we
defined }\sum_{1\leq g_{1}<g_{2}<\cdots<g_{n}\leq m}\\\text{to be an
abbreviation for}\\\sum_{\substack{\left(  g_{1},g_{2},\ldots,g_{n}\right)
\in\left\{  1,2,\ldots,m\right\}  ^{n};\\g_{1}<g_{2}<\cdots<g_{n}}}\text{)}%
}}\det\left(  \operatorname*{cols}\nolimits_{g_{1},g_{2},\ldots,g_{n}%
}A\right)  \cdot\det\left(  \operatorname*{rows}\nolimits_{g_{1},g_{2}%
,\ldots,g_{n}}B\right) \\
&  =\sum_{1\leq g_{1}<g_{2}<\cdots<g_{n}\leq m}\det\left(
\operatorname*{cols}\nolimits_{g_{1},g_{2},\ldots,g_{n}}A\right)  \cdot
\det\left(  \operatorname*{rows}\nolimits_{g_{1},g_{2},\ldots,g_{n}}B\right)
.
\end{align*}
This proves Theorem \ref{thm.cauchy-binet}.
\end{proof}
\end{verlong}

\subsection{Prelude to Laplace expansion}

Next we shall show a fact which will allow us to compute some determinants by induction:

\begin{theorem}
\label{thm.laplace.pre}Let $n$ be a positive integer. Let $A=\left(
a_{i,j}\right)  _{1\leq i\leq n,\ 1\leq j\leq n}$ be an $n\times n$-matrix.
Assume that%
\begin{equation}
a_{n,j}=0\ \ \ \ \ \ \ \ \ \ \text{for every }j\in\left\{  1,2,\ldots
,n-1\right\}  . \label{eq.thm.laplace.pre.ass}%
\end{equation}
Then, $\det A=a_{n,n}\cdot\det\left(  \left(  a_{i,j}\right)  _{1\leq i\leq
n-1,\ 1\leq j\leq n-1}\right)  $.
\end{theorem}

The assumption (\ref{eq.thm.laplace.pre.ass}) says that the last row of the
matrix $A$ consists entirely of zeroes, apart from its last entry $a_{n,n}$
(which may and may not be $0$). Theorem \ref{thm.laplace.pre} states that,
under this assumption, the determinant can be obtained by multiplying this
last entry $a_{n,n}$ with the determinant of the $\left(  n-1\right)
\times\left(  n-1\right)  $-matrix obtained by removing both the last row and
the last column from $A$. For example, for $n=3$, Theorem
\ref{thm.laplace.pre} states that%
\[
\det\left(
\begin{array}
[c]{ccc}%
a & b & c\\
d & e & f\\
0 & 0 & g
\end{array}
\right)  =g\det\left(
\begin{array}
[c]{cc}%
a & b\\
d & e
\end{array}
\right)  .
\]

Theorem \ref{thm.laplace.pre} is a particular case of \textit{Laplace
expansion}, which is a general recursive formula for the determinants that we
will encounter further below. But Theorem \ref{thm.laplace.pre} already has
noticeable applications of its own, which is why I have chosen to start with
this particular case.

The proof of Theorem \ref{thm.laplace.pre} essentially relies on the following fact:

\begin{lemma}
\label{lem.laplace.lem}Let $n$ be a positive integer. Let $\left(
a_{i,j}\right)  _{1\leq i\leq n-1,\ 1\leq j\leq n-1}$ be an $\left(
n-1\right)  \times\left(  n-1\right)  $-matrix. Then,%
\[
\sum_{\substack{\sigma\in S_{n};\\\sigma\left(  n\right)  =n}}\left(
-1\right)  ^{\sigma}\prod_{i=1}^{n-1}a_{i,\sigma\left(  i\right)  }%
=\det\left(  \left(  a_{i,j}\right)  _{1\leq i\leq n-1,\ 1\leq j\leq
n-1}\right)  .
\]

\end{lemma}

\begin{proof}
[Proof of Lemma \ref{lem.laplace.lem}.]We define a subset $T$ of $S_{n}$ by%
\[
T=\left\{  \tau\in S_{n}\ \mid\ \tau\left(  n\right)  =n\right\}  .
\]

(In other words, $T$ is the set of all $\tau\in S_{n}$ such that if we write
$\tau$ in one-line notation, then $\tau$ ends with an $n$.)

Now, we shall construct two mutually inverse maps between $S_{n-1}$ and $T$.

\begin{vershort}
For every $\sigma\in S_{n-1}$, we define a map $\widehat{\sigma}:\left\{
1,2,\ldots,n\right\}  \rightarrow\left\{  1,2,\ldots,n\right\}  $ by setting%
\[
\widehat{\sigma}\left(  i\right)  =%
\begin{cases}
\sigma\left(  i\right)  , & \text{if }i<n;\\
n, & \text{if }i=n
\end{cases}
\ \ \ \ \ \ \ \ \ \ \text{for every }i\in\left\{  1,2,\ldots,n\right\}  .
\]
\footnote{Note that if we use Definition \ref{def.perm.extend.YtX}, then this
map $\widehat{\sigma}$ is exactly the map $\sigma^{\left(  \left\{
1,2,\ldots,n-1\right\}  \rightarrow\left\{  1,2,\ldots,n\right\}  \right)  }$
.} It is straightforward to see that this map $\widehat{\sigma}$ is
well-defined and belongs to $T$. Thus, we can define a map $\Phi
:S_{n-1}\rightarrow T$ by setting%
\[
\Phi\left(  \sigma\right)  =\widehat{\sigma}\ \ \ \ \ \ \ \ \ \ \text{for
every }\sigma\in S_{n-1}.
\]

\end{vershort}

\begin{verlong}
For every $\sigma\in S_{n-1}$, we define a map $\widehat{\sigma}:\left\{
1,2,\ldots,n\right\}  \rightarrow\left\{  1,2,\ldots,n\right\}  $ by setting%
\[
\left(  \widehat{\sigma}\left(  i\right)  =%
\begin{cases}
\sigma\left(  i\right)  , & \text{if }i<n;\\
n, & \text{if }i=n
\end{cases}
\ \ \ \ \ \ \ \ \ \ \text{for every }i\in\left\{  1,2,\ldots,n\right\}
\right)  .
\]
\footnote{Note that if we use Definition \ref{def.perm.extend.YtX}, then this
map $\widehat{\sigma}$ is exactly the map $\sigma^{\left(  \left\{
1,2,\ldots,n-1\right\}  \rightarrow\left\{  1,2,\ldots,n\right\}  \right)  }$
.} This map $\widehat{\sigma}$ is well-defined\footnote{\textit{Proof.} Let
$\sigma\in S_{n-1}$. Thus, $\sigma$ is a permutation of the set $\left\{
1,2,\ldots,n-1\right\}  $ (since $S_{n-1}$ is the set of all permutations of
the set $\left\{  1,2,\ldots,n-1\right\}  $). In other words, $\sigma$ is a
bijection $\left\{  1,2,\ldots,n-1\right\}  \rightarrow\left\{  1,2,\ldots
,n-1\right\}  $.
\par
We have $n\in\left\{  1,2,\ldots,n\right\}  $ (since $n$ is positive). Now,
for every $i\in\left\{  1,2,\ldots,n\right\}  $, we have%
\begin{align*}
\widehat{\sigma}\left(  i\right)   &  =%
\begin{cases}
\sigma\left(  i\right)  , & \text{if }i<n;\\
n, & \text{if }i=n
\end{cases}
\in%
\begin{cases}
\left\{  1,2,\ldots,n\right\}  , & \text{if }i<n;\\
\left\{  1,2,\ldots,n\right\}  , & \text{if }i=n
\end{cases}
\\
&  \ \ \ \ \ \ \ \ \ \ \left(
\begin{array}
[c]{c}%
\text{since }\sigma\left(  i\right)  \in\left\{  1,2,\ldots,n-1\right\}
\subseteq\left\{  1,2,\ldots,n\right\}  \text{ in the case when }i<n\text{,}\\
\text{and since }n\in\left\{  1,2,\ldots,n\right\}  \text{ in the case when
}i=n
\end{array}
\right) \\
&  =\left\{  1,2,\ldots,n\right\}  .
\end{align*}
In other words, the map $\widehat{\sigma}$ is well-defined.} and belongs to
$T$\ \ \ \ \footnote{\textit{Proof.} Let $\sigma\in S_{n-1}$. Thus, $\sigma$
is a permutation of the set $\left\{  1,2,\ldots,n-1\right\}  $ (since
$S_{n-1}$ is the set of all permutations of the set $\left\{  1,2,\ldots
,n-1\right\}  $). In other words, $\sigma$ is a bijection $\left\{
1,2,\ldots,n-1\right\}  \rightarrow\left\{  1,2,\ldots,n-1\right\}  $. Hence,
the map $\sigma$ is injective and surjective.
\par
The definition of $\widehat{\sigma}$ yields $\widehat{\sigma}\left(  n\right)
=%
\begin{cases}
\sigma\left(  n\right)  , & \text{if }n<n;\\
n, & \text{if }n=n
\end{cases}
=n$ (since $n=n$).
\par
Let us now show that $\widehat{\sigma}$ is injective.
\par
\textit{Proof that }$\widehat{\sigma}$ \textit{is injective:} Let $p$ and $q$
be two elements of $\left\{  1,2,\ldots,n\right\}  $ such that
$\widehat{\sigma}\left(  p\right)  =\widehat{\sigma}\left(  q\right)  $. We
shall prove that $p=q$.
\par
We can WLOG assume that $p\geq q$ (since otherwise, we can simply swap $p$
with $q$). Assume this.
\par
Let us first assume that $p=n$. Then, $\widehat{\sigma}\left(  \underbrace{p}%
_{=n}\right)  =\widehat{\sigma}\left(  n\right)  =n$. Therefore, if we had
$q<n$, then we would have%
\begin{align*}
n  &  =\widehat{\sigma}\left(  p\right)  =\widehat{\sigma}\left(  q\right)  =%
\begin{cases}
\sigma\left(  q\right)  , & \text{if }q<n;\\
q, & \text{if }q=n
\end{cases}
=\sigma\left(  q\right)  \ \ \ \ \ \ \ \ \ \ \left(  \text{since }q<n\right)
\\
&  <n\ \ \ \ \ \ \ \ \ \ \left(  \text{since }\sigma\left(  q\right)
\in\left\{  1,2,\ldots,n-1\right\}  \right)  ,
\end{align*}
which is absurd. Hence, we cannot have $q<n$. Thus, we must have $q\geq n$.
Therefore, $q=n$ (since $q\in\left\{  1,2,\ldots,n\right\}  $), so that
$p=n=q$.
\par
Now, let us forget that we have assumed that $p=n$. We thus have shown that
$p=q$ in the case when $p=n$. Hence, for the rest of this proof of $p=q$, we
can WLOG assume that we don't have $p=n$. Assume this.
\par
We have $p\neq n$ (since we don't have $p=n$). Since $p\in\left\{
1,2,\ldots,n\right\}  $ and $p\neq n$, we have $p\in\left\{  1,2,\ldots
,n\right\}  \setminus\left\{  n\right\}  =\left\{  1,2,\ldots,n-1\right\}  $,
so that $p<n$. The definition of $\widehat{\sigma}$ yields $\widehat{\sigma
}\left(  p\right)  =%
\begin{cases}
\sigma\left(  p\right)  , & \text{if }p<n;\\
p, & \text{if }p=n
\end{cases}
=\sigma\left(  p\right)  $ (since $p<n$).
\par
Now, $p\geq q$, so that $q\leq p<n$. The definition of $\widehat{\sigma}$
yields $\widehat{\sigma}\left(  q\right)  =%
\begin{cases}
\sigma\left(  q\right)  , & \text{if }q<n;\\
q, & \text{if }q=n
\end{cases}
=\sigma\left(  q\right)  $ (since $q<n$). Also, since $q<n$, we have
$q\in\left\{  1,2,\ldots,n-1\right\}  $.
\par
Now, $\sigma\left(  p\right)  =\widehat{\sigma}\left(  p\right)
=\widehat{\sigma}\left(  q\right)  =\sigma\left(  q\right)  $. Hence, $p=q$
(since the map $\sigma$ is injective). This completes our proof of $p=q$.
\par
Now, let us forget that we fixed $p$ and $q$. We thus have shown that if $p$
and $q$ are two elements of $\left\{  1,2,\ldots,n\right\}  $ such that
$\widehat{\sigma}\left(  p\right)  =\widehat{\sigma}\left(  q\right)  $, then
$p=q$. In other words, the map $\widehat{\sigma}$ is injective.
\par
\textit{Proof that }$\widehat{\sigma}$ \textit{is surjective:} Let
$g\in\left\{  1,2,\ldots,n\right\}  $. We shall show that $g\in\widehat{\sigma
}\left(  \left\{  1,2,\ldots,n\right\}  \right)  $.
\par
Indeed, if $g=n$, then $g=n=\widehat{\sigma}\left(  \underbrace{n}%
_{\in\left\{  1,2,\ldots,n\right\}  }\right)  \in\widehat{\sigma}\left(
\left\{  1,2,\ldots,n\right\}  \right)  $. Hence, for the rest of this proof
of $g\in\widehat{\sigma}\left(  \left\{  1,2,\ldots,n\right\}  \right)  $, we
can WLOG assume that we don't have $g=n$. Assume this.
\par
We have $g\neq n$ (since we don't have $g=n$). Since $g\in\left\{
1,2,\ldots,n\right\}  $ and $g\neq n$, we have $g\in\left\{  1,2,\ldots
,n\right\}  \setminus\left\{  n\right\}  =\left\{  1,2,\ldots,n-1\right\}
=\sigma\left(  \left\{  1,2,\ldots,n-1\right\}  \right)  $ (since the map
$\sigma$ is surjective). In other words, there exists an $h\in\left\{
1,2,\ldots,n-1\right\}  $ such that $g=\sigma\left(  h\right)  $. Consider
this $h$. We have $h\in\left\{  1,2,\ldots,n-1\right\}  \subseteq\left\{
1,2,\ldots,n\right\}  $ and $h<n$ (since $h\in\left\{  1,2,\ldots,n-1\right\}
$). The definition of $\widehat{\sigma}$ yields $\widehat{\sigma}\left(
h\right)  =%
\begin{cases}
\sigma\left(  h\right)  , & \text{if }h<n;\\
h, & \text{if }h=n
\end{cases}
=\sigma\left(  h\right)  $ (since $h<n$). Compared with $g=\sigma\left(
h\right)  $, this yields $g=\widehat{\sigma}\left(  \underbrace{h}%
_{\in\left\{  1,2,\ldots,n\right\}  }\right)  \in\widehat{\sigma}\left(
\left\{  1,2,\ldots,n\right\}  \right)  $. Thus, $g\in\widehat{\sigma}\left(
\left\{  1,2,\ldots,n\right\}  \right)  $ is proven.
\par
Now, let us forget that we fixed $g$. We thus have shown that $g\in
\widehat{\sigma}\left(  \left\{  1,2,\ldots,n\right\}  \right)  $ for every
$g\in\left\{  1,2,\ldots,n\right\}  $. In other words, $\left\{
1,2,\ldots,n\right\}  \subseteq\widehat{\sigma}\left(  \left\{  1,2,\ldots
,n\right\}  \right)  $. In other words, the map $\widehat{\sigma}$ is
surjective.
\par
We now know that the map $\widehat{\sigma}$ is both injective and surjective.
In other words, $\widehat{\sigma}$ is bijective. Hence, $\widehat{\sigma}$ is
a bijective map $\left\{  1,2,\ldots,n\right\}  \rightarrow\left\{
1,2,\ldots,n\right\}  $. In other words, $\widehat{\sigma}$ is a permutation
of $\left\{  1,2,\ldots,n\right\}  $. In other words, $\widehat{\sigma}\in
S_{n}$ (since $S_{n}$ is the set of all permutations of the set $\left\{
1,2,\ldots,n\right\}  $).
\par
Thus, $\widehat{\sigma}$ is an element of $S_{n}$ and satisfies
$\widehat{\sigma}\left(  n\right)  =n$. In other words, $\widehat{\sigma}$ is
an element $\tau$ of $S_{n}$ satisfying $\tau\left(  n\right)  =n$. In other
words,%
\[
\widehat{\sigma}\in\left\{  \tau\in S_{n}\ \mid\ \tau\left(  n\right)
=n\right\}  =T,
\]
qed.}. Thus, we can define a map $\Phi:S_{n-1}\rightarrow T$ by setting%
\[
\left(  \Phi\left(  \sigma\right)  =\widehat{\sigma}%
\ \ \ \ \ \ \ \ \ \ \text{for every }\sigma\in S_{n-1}\right)  .
\]

\end{verlong}

Loosely speaking, for every $\sigma\in S_{n-1}$, the permutation $\Phi\left(
\sigma\right)  =\widehat{\sigma}\in T$ is obtained by writing $\sigma$ in
one-line notation and appending $n$ on its right. For example, if $n=4$ and if
$\sigma\in S_{3}$ is the permutation that is written as $\left(  2,3,1\right)
$ in one-line notation, then $\Phi\left(  \sigma\right)  =\widehat{\sigma}$ is
the permutation that is written as $\left(  2,3,1,4\right)  $ in one-line notation.

\begin{vershort}
On the other hand, for every $\gamma\in T$, we define a map $\overline{\gamma
}:\left\{  1,2,\ldots,n-1\right\}  \rightarrow\left\{  1,2,\ldots,n-1\right\}
$ by setting
\[
\overline{\gamma}\left(  i\right)  =\gamma\left(  i\right)
\ \ \ \ \ \ \ \ \ \ \text{for every }i\in\left\{  1,2,\ldots,n-1\right\}  .
\]
It is straightforward to see that this map $\overline{\gamma}$ is well-defined
and belongs to $S_{n-1}$. Hence, we can define a map $\Psi:T\rightarrow
S_{n-1}$ by setting%
\[
\Psi\left(  \gamma\right)  =\overline{\gamma}\ \ \ \ \ \ \ \ \ \ \text{for
every }\gamma\in T.
\]

\end{vershort}

\begin{verlong}
On the other hand, for every $\gamma\in T$, we define a map $\overline{\gamma
}:\left\{  1,2,\ldots,n-1\right\}  \rightarrow\left\{  1,2,\ldots,n-1\right\}
$ by setting
\[
\left(  \overline{\gamma}\left(  i\right)  =\gamma\left(  i\right)
\ \ \ \ \ \ \ \ \ \ \text{for every }i\in\left\{  1,2,\ldots,n-1\right\}
\right)  .
\]
This map $\overline{\gamma}$ is well-defined\footnote{\textit{Proof.} Let
$\gamma\in T$. Thus, $\gamma\in T=\left\{  \tau\in S_{n}\ \mid\ \tau\left(
n\right)  =n\right\}  $. In other words, $\gamma$ is an element $\tau$ of
$S_{n}$ satisfying $\tau\left(  n\right)  =n$. In other words, $\gamma$ is an
element of $S_{n}$ and satisfies $\gamma\left(  n\right)  =n$.
\par
We have $\gamma\in S_{n}$. In other words, $\gamma$ is a permutation of the
set $\left\{  1,2,\ldots,n\right\}  $ (since $S_{n}$ is the set of all
permutations of the set $\left\{  1,2,\ldots,n\right\}  $). In other words,
$\gamma$ is a bijection $\left\{  1,2,\ldots,n\right\}  \rightarrow\left\{
1,2,\ldots,n\right\}  $. Hence, the map $\gamma$ is surjective and injective.
\par
Let $i\in\left\{  1,2,\ldots,n-1\right\}  $. Then, $i<n$ and thus $i\neq n$.
Also, $i\in\left\{  1,2,\ldots,n-1\right\}  \subseteq\left\{  1,2,\ldots
,n\right\}  $. Thus, $\gamma\left(  i\right)  $ is well-defined and belongs to
$\left\{  1,2,\ldots,n\right\}  $. Also, $i\neq n$. Hence, $\gamma\left(
i\right)  \neq\gamma\left(  n\right)  $ (since the map $\gamma$ is injective),
so that $\gamma\left(  i\right)  \neq\gamma\left(  n\right)  =n$. Combined
with $\gamma\left(  i\right)  \in\left\{  1,2,\ldots,n\right\}  $, this shows
that $\gamma\left(  i\right)  \in\left\{  1,2,\ldots,n\right\}  \setminus
\left\{  n\right\}  =\left\{  1,2,\ldots,n-1\right\}  $.
\par
Now, let us forget that we fixed $i$. We thus have shown that $\gamma\left(
i\right)  \in\left\{  1,2,\ldots,n-1\right\}  $ for every $i\in\left\{
1,2,\ldots,n-1\right\}  $. In other words, the map $\overline{\gamma}$ is
well-defined. Qed.} and belongs to $S_{n-1}$\ \ \ \ \footnote{\textit{Proof.}
Let $\gamma\in T$. Thus, $\gamma\in T=\left\{  \tau\in S_{n}\ \mid
\ \tau\left(  n\right)  =n\right\}  $. In other words, $\gamma$ is an element
$\tau$ of $S_{n}$ satisfying $\tau\left(  n\right)  =n$. In other words,
$\gamma$ is an element of $S_{n}$ and satisfies $\gamma\left(  n\right)  =n$.
\par
We have $\gamma\in S_{n}$. In other words, $\gamma$ is a permutation of the
set $\left\{  1,2,\ldots,n\right\}  $ (since $S_{n}$ is the set of all
permutations of the set $\left\{  1,2,\ldots,n\right\}  $). In other words,
$\gamma$ is a bijection $\left\{  1,2,\ldots,n\right\}  \rightarrow\left\{
1,2,\ldots,n\right\}  $. Hence, the map $\gamma$ is surjective and injective.
\par
Let us now show that $\overline{\gamma}$ is injective.
\par
\textit{Proof that }$\overline{\gamma}$ \textit{is injective:} Let $p$ and $q$
be two elements of $\left\{  1,2,\ldots,n-1\right\}  $ such that
$\overline{\gamma}\left(  p\right)  =\overline{\gamma}\left(  q\right)  $. We
shall prove that $p=q$.
\par
We have $\overline{\gamma}\left(  p\right)  =\gamma\left(  p\right)  $ (by the
definition of $\overline{\gamma}$) and $\overline{\gamma}\left(  q\right)
=\gamma\left(  q\right)  $ (by the definition of $\overline{\gamma}$). Thus,
$\gamma\left(  p\right)  =\overline{\gamma}\left(  p\right)  =\overline
{\gamma}\left(  q\right)  =\gamma\left(  q\right)  $. Therefore, $p=q$ (since
the map $\gamma$ is injective).
\par
Now, let us forget that we fixed $p$ and $q$. We thus have shown that if $p$
and $q$ are two elements of $\left\{  1,2,\ldots,n\right\}  $ such that
$\overline{\gamma}\left(  p\right)  =\overline{\gamma}\left(  q\right)  $,
then $p=q$. In other words, the map $\overline{\gamma}$ is injective.
\par
\textit{Proof that }$\overline{\gamma}$ \textit{is surjective:} Let
$g\in\left\{  1,2,\ldots,n-1\right\}  $. We shall show that $g\in
\overline{\gamma}\left(  \left\{  1,2,\ldots,n-1\right\}  \right)  $.
\par
Indeed, $g\in\left\{  1,2,\ldots,n-1\right\}  \subseteq\left\{  1,2,\ldots
,n\right\}  =\gamma\left(  \left\{  1,2,\ldots,n\right\}  \right)  $ (since
the map $\gamma$ is surjective). Hence, there exists an $h\in\left\{
1,2,\ldots,n\right\}  $ such that $g=\gamma\left(  h\right)  $. Consider this
$h$.
\par
We have $g\in\left\{  1,2,\ldots,n-1\right\}  $, thus $g<n$, thus $g\neq n$.
\par
If we had $h=n$, then we would have $g=\gamma\left(  \underbrace{h}%
_{=n}\right)  =\gamma\left(  n\right)  =n$, which would contradict $g\neq n$.
Hence, we cannot have $h=n$. We thus have $h\neq n$. Combined with
$h\in\left\{  1,2,\ldots,n\right\}  $, this shows that $h\in\left\{
1,2,\ldots,n\right\}  \setminus\left\{  n\right\}  =\left\{  1,2,\ldots
,n-1\right\}  $. Thus, $\overline{\gamma}\left(  h\right)  $ is well-defined.
The definition of $\overline{\gamma}$ shows that $\overline{\gamma}\left(
h\right)  =\gamma\left(  h\right)  =g$.
\par
Hence, $g=\overline{\gamma}\left(  \underbrace{h}_{\in\left\{  1,2,\ldots
,n-1\right\}  }\right)  \in\overline{\gamma}\left(  \left\{  1,2,\ldots
,n-1\right\}  \right)  $.
\par
Now, let us forget that we fixed $g$. We thus have shown that $g\in
\overline{\gamma}\left(  \left\{  1,2,\ldots,n-1\right\}  \right)  $ for every
$g\in\left\{  1,2,\ldots,n-1\right\}  $. In other words, $\left\{
1,2,\ldots,n-1\right\}  \subseteq\overline{\gamma}\left(  \left\{
1,2,\ldots,n-1\right\}  \right)  $. In other words, the map $\overline{\gamma
}$ is surjective.
\par
We now know that the map $\overline{\gamma}$ is both injective and surjective.
In other words, $\overline{\gamma}$ is bijective. Hence, $\overline{\gamma}$
is a bijective map $\left\{  1,2,\ldots,n-1\right\}  \rightarrow\left\{
1,2,\ldots,n-1\right\}  $. In other words, $\overline{\gamma}$ is a
permutation of $\left\{  1,2,\ldots,n-1\right\}  $. In other words,
$\overline{\gamma}\in S_{n-1}$ (since $S_{n-1}$ is the set of all permutations
of the set $\left\{  1,2,\ldots,n-1\right\}  $), qed.}. Hence, we can define a
map $\Psi:T\rightarrow S_{n-1}$ by setting%
\[
\left(  \Psi\left(  \gamma\right)  =\overline{\gamma}%
\ \ \ \ \ \ \ \ \ \ \text{for every }\gamma\in T\right)  .
\]

\end{verlong}

Loosely speaking, for every $\gamma\in T$, the permutation $\Psi\left(
\gamma\right)  =\overline{\gamma}\in S_{n-1}$ is obtained by writing $\gamma$
in one-line notation and removing the $n$ (which is the rightmost entry in the
one-line notation, because $\gamma\left(  n\right)  =n$). For example, if
$n=4$ and if $\gamma\in S_{4}$ is the permutation that is written as $\left(
2,3,1,4\right)  $ in one-line notation, then $\Psi\left(  \gamma\right)
=\overline{\gamma}$ is the permutation that is written as $\left(
2,3,1\right)  $ in one-line notation.

\begin{vershort}
The maps $\Phi$ and $\Psi$ are mutually inverse.\footnote{This should be clear
enough from the descriptions we gave using one-line notation. A formal proof
is straightforward.} Thus, the map $\Phi$ is a bijection.
\end{vershort}

\begin{verlong}
We have $\Phi\circ\Psi=\operatorname*{id}$\ \ \ \ \footnote{\textit{Proof.}
Let $\gamma\in T$. Then, $\gamma\in T=\left\{  \tau\in S_{n}\ \mid
\ \tau\left(  n\right)  =n\right\}  $. In other words, $\gamma$ is an element
$\tau$ of $S_{n}$ satisfying $\tau\left(  n\right)  =n$. In other words,
$\gamma$ is an element of $S_{n}$ and satisfies $\gamma\left(  n\right)  =n$.
\par
Let $\sigma=\overline{\gamma}$. Then, the definition of $\Phi$ yields
$\Phi\left(  \sigma\right)  =\widehat{\sigma}$. Also, the definition of $\Psi$
yields $\Psi\left(  \gamma\right)  =\overline{\gamma}=\sigma$. Now, $\left(
\Phi\circ\Psi\right)  \left(  \gamma\right)  =\Phi\left(  \underbrace{\Psi
\left(  \gamma\right)  }_{=\sigma}\right)  =\Phi\left(  \sigma\right)
=\widehat{\sigma}$.
\par
Both $\widehat{\sigma}$ and $\gamma$ are elements of $T$, therefore elements
of $S_{n}$ (since $T\subseteq S_{n}$), therefore permutations of the set
$\left\{  1,2,\ldots,n\right\}  $ (since $S_{n}$ is the set of all
permutations of the set $\left\{  1,2,\ldots,n\right\}  $), therefore
bijective maps $\left\{  1,2,\ldots,n\right\}  \rightarrow\left\{
1,2,\ldots,n\right\}  $.
\par
Let $i\in\left\{  1,2,\ldots,n\right\}  $. We shall show that $\widehat{\sigma
}\left(  i\right)  =\gamma\left(  i\right)  $.
\par
If $i=n$, then%
\begin{align*}
\widehat{\sigma}\left(  i\right)   &  =%
\begin{cases}
\sigma\left(  i\right)  , & \text{if }i<n;\\
n, & \text{if }i=n
\end{cases}
\ \ \ \ \ \ \ \ \ \ \left(  \text{by the definition of }\widehat{\sigma
}\right) \\
&  =n\ \ \ \ \ \ \ \ \ \ \left(  \text{since }i=n\right) \\
&  =\gamma\left(  \underbrace{n}_{=i}\right)  \ \ \ \ \ \ \ \ \ \ \left(
\text{since }\gamma\left(  n\right)  =n\right) \\
&  =\gamma\left(  i\right)  .
\end{align*}
Hence, $\widehat{\sigma}\left(  i\right)  =\gamma\left(  i\right)  $ is proven
in the case when $i=n$. Therefore, for the rest of the proof of
$\widehat{\sigma}\left(  i\right)  =\gamma\left(  i\right)  $, we can WLOG
assume that we don't have $i=n$. Assume this.
\par
We have $i\neq n$ (since we don't have $i=n$). Combined with $i\in\left\{
1,2,\ldots,n\right\}  $, this shows that $i\in\left\{  1,2,\ldots,n\right\}
\setminus\left\{  n\right\}  =\left\{  1,2,\ldots,n-1\right\}  $. Hence,
$i<n$. Now,%
\begin{align*}
\widehat{\sigma}\left(  i\right)   &  =%
\begin{cases}
\sigma\left(  i\right)  , & \text{if }i<n;\\
n, & \text{if }i=n
\end{cases}
\ \ \ \ \ \ \ \ \ \ \left(  \text{by the definition of }\widehat{\sigma
}\right) \\
&  =\underbrace{\sigma}_{=\overline{\gamma}}\left(  i\right)
\ \ \ \ \ \ \ \ \ \ \left(  \text{since }i<n\right) \\
&  =\overline{\gamma}\left(  i\right)  =\gamma\left(  i\right)
\ \ \ \ \ \ \ \ \ \ \left(  \text{by the definition of }\overline{\gamma
}\right)  .
\end{align*}
Thus, $\widehat{\sigma}\left(  i\right)  =\gamma\left(  i\right)  $ is proven.
\par
Let us now forget that we fixed $i$. We thus have proven that $\widehat{\sigma
}\left(  i\right)  =\gamma\left(  i\right)  $ for every $i\in\left\{
1,2,\ldots,n\right\}  $. In other words, $\widehat{\sigma}=\gamma$ (since both
$\widehat{\sigma}$ and $\gamma$ are maps $\left\{  1,2,\ldots,n\right\}
\rightarrow\left\{  1,2,\ldots,n\right\}  $). Hence, $\left(  \Phi\circ
\Psi\right)  \left(  \gamma\right)  =\widehat{\sigma}=\gamma
=\operatorname*{id}\left(  \gamma\right)  $.
\par
Let us now forget that we fixed $\gamma$. We have thus proven that $\left(
\Phi\circ\Psi\right)  \left(  \gamma\right)  =\operatorname*{id}\left(
\gamma\right)  $ for every $\gamma\in T$. In other words, $\Phi\circ
\Psi=\operatorname*{id}$, qed.} and $\Psi\circ\Phi=\operatorname*{id}%
$\ \ \ \ \footnote{\textit{Proof.} Let $\sigma\in S_{n-1}$. Let $\varepsilon
=\widehat{\sigma}$. Then, the definition of $\Psi$ yields $\Phi\left(
\sigma\right)  =\widehat{\sigma}=\varepsilon$. Also, the definition of $\Psi$
yields $\Psi\left(  \varepsilon\right)  =\overline{\varepsilon}$. Now,
$\left(  \Psi\circ\Phi\right)  \left(  \sigma\right)  =\Psi\left(
\underbrace{\Phi\left(  \sigma\right)  }_{=\varepsilon}\right)  =\Psi\left(
\varepsilon\right)  =\overline{\varepsilon}$.
\par
Both $\overline{\varepsilon}$ and $\sigma$ are elements of $S_{n-1}$,
therefore permutations of the set $\left\{  1,2,\ldots,n-1\right\}  $ (since
$S_{n-1}$ is the set of all permutations of the set $\left\{  1,2,\ldots
,n-1\right\}  $), therefore bijective maps $\left\{  1,2,\ldots,n-1\right\}
\rightarrow\left\{  1,2,\ldots,n-1\right\}  $.
\par
Let $i\in\left\{  1,2,\ldots,n-1\right\}  $. Thus, $i\leq n-1<n$. Now,%
\begin{align*}
\overline{\varepsilon}\left(  i\right)   &  =\underbrace{\varepsilon
}_{=\widehat{\sigma}}\left(  i\right)  \ \ \ \ \ \ \ \ \ \ \left(  \text{by
the definition of }\overline{\varepsilon}\right) \\
&  =\widehat{\sigma}\left(  i\right)  =%
\begin{cases}
\sigma\left(  i\right)  , & \text{if }i<n;\\
n, & \text{if }i=n
\end{cases}
\ \ \ \ \ \ \ \ \ \ \left(  \text{by the definition of }\widehat{\sigma
}\right) \\
&  =\sigma\left(  i\right)  \ \ \ \ \ \ \ \ \ \ \left(  \text{since
}i<n\right)  .
\end{align*}
\par
Let us now forget that we fixed $i$. We thus have proven that $\overline
{\varepsilon}\left(  i\right)  =\sigma\left(  i\right)  $ for every
$i\in\left\{  1,2,\ldots,n-1\right\}  $. In other words, $\overline
{\varepsilon}=\sigma$ (since both $\overline{\varepsilon}$ and $\sigma$ are
maps $\left\{  1,2,\ldots,n-1\right\}  \rightarrow\left\{  1,2,\ldots
,n-1\right\}  $). Hence, $\left(  \Psi\circ\Phi\right)  \left(  \sigma\right)
=\overline{\varepsilon}=\sigma=\operatorname*{id}\left(  \sigma\right)  $.
\par
Let us now forget that we fixed $\sigma$. We have thus proven that $\left(
\Psi\circ\Phi\right)  \left(  \sigma\right)  =\operatorname*{id}\left(
\sigma\right)  $ for every $\sigma\in S_{n-1}$. In other words, $\Psi\circ
\Phi=\operatorname*{id}$, qed.}. These two equalities show that the maps
$\Phi$ and $\Psi$ are mutually inverse. Hence, the map $\Phi$ is a bijection.
\end{verlong}

\begin{vershort}
It is fairly easy to see that every $\sigma\in S_{n-1}$ satisfies%
\begin{equation}
\left(  -1\right)  ^{\widehat{\sigma}}=\left(  -1\right)  ^{\sigma}
\label{pf.lem.laplace.lem.short.-1}%
\end{equation}
\footnote{\textit{Proof of (\ref{pf.lem.laplace.lem.short.-1}):} Let
$\sigma\in S_{n-1}$. We want to prove that $\left(  -1\right)
^{\widehat{\sigma}}=\left(  -1\right)  ^{\sigma}$. It is clearly sufficient to
show that $\ell\left(  \widehat{\sigma}\right)  =\ell\left(  \sigma\right)  $
(because $\left(  -1\right)  ^{\widehat{\sigma}}=\left(  -1\right)
^{\ell\left(  \widehat{\sigma}\right)  }$ and $\left(  -1\right)  ^{\sigma
}=\left(  -1\right)  ^{\ell\left(  \sigma\right)  }$). In order to do so, it
is sufficient to show that the inversions of $\widehat{\sigma}$ are precisely
the inversions of $\sigma$ (because $\ell\left(  \widehat{\sigma}\right)  $ is
the number of inversions of $\widehat{\sigma}$, whereas $\ell\left(
\sigma\right)  $ is the number of inversions of $\sigma$).
\par
If $\left(  i,j\right)  $ is an inversion of $\sigma$, then $\left(
i,j\right)  $ is an inversion of $\widehat{\sigma}$ (because if $\left(
i,j\right)  $ is an inversion of $\sigma$, then both $i$ and $j$ are $<n$, and
thus the definition of $\widehat{\sigma}$ yields $\widehat{\sigma}\left(
i\right)  =\sigma\left(  i\right)  $ and $\widehat{\sigma}\left(  j\right)
=\sigma\left(  j\right)  $). In other words, every inversion of $\sigma$ is an
inversion of $\widehat{\sigma}$.
\par
On the other hand, let $\left(  u,v\right)  $ be an inversion of
$\widehat{\sigma}$. We shall prove that $\left(  u,v\right)  $ is an inversion
of $\sigma$.
\par
Indeed, $\left(  u,v\right)  $ is an inversion of $\widehat{\sigma}$. In other
words, $\left(  u,v\right)  $ is a pair of integers satisfying $1\leq u<v\leq
n$ and $\widehat{\sigma}\left(  u\right)  >\widehat{\sigma}\left(  v\right)
$.
\par
If we had $v=n$, then we would have $\widehat{\sigma}\left(  u\right)
>\widehat{\sigma}\left(  \underbrace{v}_{=n}\right)  =\widehat{\sigma}\left(
n\right)  =n$ (by the definition of $\widehat{\sigma}$), which would
contradict $\widehat{\sigma}\left(  u\right)  \in\left\{  1,2,\ldots
,n\right\}  $. Thus, we cannot have $v=n$. We therefore have $v<n$, so that
$v\leq n-1$. Now, $1\leq u<v\leq n-1$. Thus, both $\sigma\left(  u\right)  $
and $\sigma\left(  v\right)  $ are well-defined. The definition of
$\widehat{\sigma}$ yields $\widehat{\sigma}\left(  u\right)  =\sigma\left(
u\right)  $ (since $u\leq n-1<n$) and $\widehat{\sigma}\left(  v\right)
=\sigma\left(  v\right)  $ (since $v\leq n-1<n$), so that $\sigma\left(
u\right)  =\widehat{\sigma}\left(  u\right)  >\widehat{\sigma}\left(
v\right)  =\sigma\left(  v\right)  $. Thus, $\left(  u,v\right)  $ is a pair
of integers satisfying $1\leq u<v\leq n-1$ and $\sigma\left(  u\right)
>\sigma\left(  v\right)  $. In other words, $\left(  u,v\right)  $ is an
inversion of $\sigma$.
\par
We thus have shown that every inversion of $\widehat{\sigma}$ is an inversion
of $\sigma$. Combining this with the fact that every inversion of $\sigma$ is
an inversion of $\widehat{\sigma}$, we thus conclude that the inversions of
$\widehat{\sigma}$ are precisely the inversions of $\sigma$. As we have
already said, this finishes the proof of (\ref{pf.lem.laplace.lem.short.-1}).}
and%
\begin{equation}
\prod_{i=1}^{n-1}a_{i,\widehat{\sigma}\left(  i\right)  }=\prod_{i=1}%
^{n-1}a_{i,\sigma\left(  i\right)  } \label{pf.lem.laplace.lem.short.prod}%
\end{equation}
\footnote{\textit{Proof of (\ref{pf.lem.laplace.lem.short.prod}):} Let
$\sigma\in S_{n-1}$. The definition of $\widehat{\sigma}$ yields
$\widehat{\sigma}\left(  i\right)  =\sigma\left(  i\right)  $ for every
$i\in\left\{  1,2,\ldots,n-1\right\}  $. Thus, $a_{i,\widehat{\sigma}\left(
i\right)  }=a_{i,\sigma\left(  i\right)  }$ for every $i\in\left\{
1,2,\ldots,n-1\right\}  $. Hence, $\prod_{i=1}^{n-1}%
\underbrace{a_{i,\widehat{\sigma}\left(  i\right)  }}_{=a_{i,\sigma\left(
i\right)  }}=\prod_{i=1}^{n-1}a_{i,\sigma\left(  i\right)  }$, qed.}.
\end{vershort}

\begin{verlong}
For every $m\in\mathbb{N}$ and every $\sigma\in S_{m}$, let
$\operatorname*{Inv}\left(  \sigma\right)  $ be the set of all inversions of
the permutation $\sigma$. Thus, for every $m\in\mathbb{N}$ and $\sigma\in
S_{m}$, we have%
\begin{align}
\ell\left(  \sigma\right)   &  =\left(  \text{the number of inversions of
}\sigma\right)  \ \ \ \ \ \ \ \ \ \ \left(  \text{by the definition of }%
\ell\left(  \sigma\right)  \right) \nonumber\\
&  =\left\vert \underbrace{\left(  \text{the set of all inversions of }%
\sigma\right)  }_{\substack{=\operatorname*{Inv}\left(  \sigma\right)
\\\text{(since }\operatorname*{Inv}\left(  \sigma\right)  \text{ is the set of
all inversions of }\sigma\text{)}}}\right\vert \nonumber\\
&  =\left\vert \operatorname*{Inv}\left(  \sigma\right)  \right\vert .
\label{pf.lem.laplace.lem.l}%
\end{align}

\end{verlong}

\begin{verlong}
Moreover, every $\sigma\in S_{n-1}$ satisfies
\[
\operatorname*{Inv}\left(  \widehat{\sigma}\right)  \subseteq
\operatorname*{Inv}\left(  \sigma\right)
\]
\footnote{\textit{Proof.} Let $\sigma\in S_{n-1}$. Let $c\in
\operatorname*{Inv}\left(  \widehat{\sigma}\right)  $. Thus, $c$ is an
inversion of $\widehat{\sigma}$ (since $\operatorname*{Inv}\left(
\widehat{\sigma}\right)  $ is the set of all inversions of $\widehat{\sigma}%
$). In other words, $c$ is a pair $\left(  i,j\right)  $ of integers
satisfying $1\leq i<j\leq n$ and $\widehat{\sigma}\left(  i\right)
>\widehat{\sigma}\left(  j\right)  $ (by the definition of an
\textquotedblleft inversion\textquotedblright).
\par
We have $\widehat{\sigma}\left(  i\right)  \in\left\{  1,2,\ldots,n\right\}
$, thus $\widehat{\sigma}\left(  i\right)  \leq n$. From $\widehat{\sigma
}\left(  i\right)  >\widehat{\sigma}\left(  j\right)  $, we obtain
$\widehat{\sigma}\left(  j\right)  <\widehat{\sigma}\left(  i\right)  \leq n$
and thus $\widehat{\sigma}\left(  j\right)  \leq n-1$ (since $\widehat{\sigma
}\left(  j\right)  $ and $n$ are integers).
\par
The definition of $\widehat{\sigma}$ yields $\widehat{\sigma}\left(  n\right)
=%
\begin{cases}
\sigma\left(  n\right)  , & \text{if }n<n;\\
n, & \text{if }n=n
\end{cases}
=n$ (since $n=n$). Hence, $\widehat{\sigma}\left(  n\right)  =n\neq
\widehat{\sigma}\left(  j\right)  $ (since $\widehat{\sigma}\left(  j\right)
<\widehat{\sigma}\left(  n\right)  $), so that $n\neq j$. Hence, $j\neq n$, so
that $j<n$ (since $j\leq n$). Hence, $1\leq i<j\leq n-1$.
\par
The definition of $\widehat{\sigma}$ yields $\widehat{\sigma}\left(  i\right)
=%
\begin{cases}
\sigma\left(  i\right)  , & \text{if }i<n;\\
n, & \text{if }i=n
\end{cases}
=\sigma\left(  i\right)  $ (since $i<n$). The definition of $\widehat{\sigma}$
yields $\widehat{\sigma}\left(  j\right)  =%
\begin{cases}
\sigma\left(  j\right)  , & \text{if }j<n;\\
n, & \text{if }j=n
\end{cases}
=\sigma\left(  j\right)  $ (since $j<n$). Now, $\sigma\left(  i\right)
=\widehat{\sigma}\left(  i\right)  >\widehat{\sigma}\left(  j\right)
=\sigma\left(  j\right)  $.
\par
Now, we know that $\left(  i,j\right)  $ is a pair of integers satisfying
$1\leq i<j\leq n-1$ and $\sigma\left(  i\right)  >\sigma\left(  j\right)  $.
In other words, $\left(  i,j\right)  $ is an inversion of $\sigma$ (by the
definition of an \textquotedblleft inversion\textquotedblright). In other
words, $\left(  i,j\right)  \in\operatorname*{Inv}\left(  \sigma\right)  $
(since $\operatorname*{Inv}\left(  \sigma\right)  $ is the set of all
inversions of $\sigma$). Thus, $c=\left(  i,j\right)  \in\operatorname*{Inv}%
\left(  \sigma\right)  $.
\par
Let us now forget that we fixed $c$. We thus have proven that $c\in
\operatorname*{Inv}\left(  \sigma\right)  $ for every $c\in\operatorname*{Inv}%
\left(  \widehat{\sigma}\right)  $. In other words, $\operatorname*{Inv}%
\left(  \widehat{\sigma}\right)  \subseteq\operatorname*{Inv}\left(
\sigma\right)  $, qed.} and%
\[
\operatorname*{Inv}\left(  \sigma\right)  \subseteq\operatorname*{Inv}\left(
\widehat{\sigma}\right)
\]
\footnote{\textit{Proof.} Let $\sigma\in S_{n-1}$. Let $c\in
\operatorname*{Inv}\left(  \sigma\right)  $. Thus, $c$ is an inversion of
$\sigma$ (since $\operatorname*{Inv}\left(  \sigma\right)  $ is the set of all
inversions of $\sigma$). In other words, $c$ is a pair $\left(  i,j\right)  $
of integers satisfying $1\leq i<j\leq n-1$ and $\sigma\left(  i\right)
>\sigma\left(  j\right)  $ (by the definition of an \textquotedblleft
inversion\textquotedblright).
\par
We have $j\leq n-1<n$ and thus $i<j<n$. Since $1\leq i\leq n-1\leq n$, we have
$i\in\left\{  1,2,\ldots,n\right\}  $. Since $1\leq j\leq n-1\leq n$, we have
$j\in\left\{  1,2,\ldots,n\right\}  $.
\par
The definition of $\widehat{\sigma}$ yields $\widehat{\sigma}\left(  i\right)
=%
\begin{cases}
\sigma\left(  i\right)  , & \text{if }i<n;\\
n, & \text{if }i=n
\end{cases}
=\sigma\left(  i\right)  $ (since $i<n$). The definition of $\widehat{\sigma}$
yields $\widehat{\sigma}\left(  j\right)  =%
\begin{cases}
\sigma\left(  j\right)  , & \text{if }j<n;\\
n, & \text{if }j=n
\end{cases}
=\sigma\left(  j\right)  $ (since $j<n$). Now, $\widehat{\sigma}\left(
i\right)  =\sigma\left(  i\right)  >\sigma\left(  j\right)  =\widehat{\sigma
}\left(  j\right)  $.
\par
Now, we know that $\left(  i,j\right)  $ is a pair of integers satisfying
$1\leq i<j\leq n$ and $\widehat{\sigma}\left(  i\right)  >\widehat{\sigma
}\left(  j\right)  $. In other words, $\left(  i,j\right)  $ is an inversion
of $\widehat{\sigma}$ (by the definition of an \textquotedblleft
inversion\textquotedblright). In other words, $\left(  i,j\right)
\in\operatorname*{Inv}\left(  \widehat{\sigma}\right)  $ (since
$\operatorname*{Inv}\left(  \widehat{\sigma}\right)  $ is the set of all
inversions of $\widehat{\sigma}$). Thus, $c=\left(  i,j\right)  \in
\operatorname*{Inv}\left(  \widehat{\sigma}\right)  $.
\par
Let us now forget that we fixed $c$. We thus have proven that $c\in
\operatorname*{Inv}\left(  \widehat{\sigma}\right)  $ for every $c\in
\operatorname*{Inv}\left(  \sigma\right)  $. In other words,
$\operatorname*{Inv}\left(  \sigma\right)  \subseteq\operatorname*{Inv}\left(
\widehat{\sigma}\right)  $, qed.} and%
\[
\operatorname*{Inv}\left(  \sigma\right)  =\operatorname*{Inv}\left(
\widehat{\sigma}\right)
\]
\footnote{\textit{Proof.} Let $\sigma\in S_{n-1}$. We have proven the two
relations $\operatorname*{Inv}\left(  \sigma\right)  \subseteq
\operatorname*{Inv}\left(  \widehat{\sigma}\right)  $ and $\operatorname*{Inv}%
\left(  \widehat{\sigma}\right)  \subseteq\operatorname*{Inv}\left(
\sigma\right)  $. Combining these two relations, we obtain
$\operatorname*{Inv}\left(  \sigma\right)  =\operatorname*{Inv}\left(
\widehat{\sigma}\right)  $. Qed.} and%
\[
\ell\left(  \sigma\right)  =\ell\left(  \widehat{\sigma}\right)
\]
\footnote{\textit{Proof.} Let $\sigma\in S_{n-1}$. We have proven that
$\operatorname*{Inv}\left(  \sigma\right)  =\operatorname*{Inv}\left(
\widehat{\sigma}\right)  $. Now, (\ref{pf.lem.laplace.lem.l}) (applied to
$m=n-1$) yields $\ell\left(  \sigma\right)  =\left\vert
\underbrace{\operatorname*{Inv}\left(  \sigma\right)  }_{=\operatorname*{Inv}%
\left(  \widehat{\sigma}\right)  }\right\vert =\left\vert \operatorname*{Inv}%
\left(  \widehat{\sigma}\right)  \right\vert $. On the other hand,
(\ref{pf.lem.laplace.lem.l}) (applied to $n-1$ and $\widehat{\sigma}$ instead
of $m$ and $\sigma$) yields $\ell\left(  \widehat{\sigma}\right)  =\left\vert
\operatorname*{Inv}\left(  \widehat{\sigma}\right)  \right\vert $. Comparing
this with $\ell\left(  \sigma\right)  =\left\vert \operatorname*{Inv}\left(
\widehat{\sigma}\right)  \right\vert $, we obtain $\ell\left(  \sigma\right)
=\ell\left(  \widehat{\sigma}\right)  $, qed.} and%
\begin{equation}
\left(  -1\right)  ^{\widehat{\sigma}}=\left(  -1\right)  ^{\sigma}
\label{pf.lem.laplace.lem.-1}%
\end{equation}
\footnote{\textit{Proof.} Let $\sigma\in S_{n-1}$. We have proven that
$\ell\left(  \sigma\right)  =\ell\left(  \widehat{\sigma}\right)  $. But the
definition of $\left(  -1\right)  ^{\sigma}$ yields $\left(  -1\right)
^{\sigma}=\left(  -1\right)  ^{\ell\left(  \sigma\right)  }$. Also, the
definition of $\left(  -1\right)  ^{\widehat{\sigma}}$ yields $\left(
-1\right)  ^{\widehat{\sigma}}=\left(  -1\right)  ^{\ell\left(
\widehat{\sigma}\right)  }$. Thus,%
\begin{align*}
\left(  -1\right)  ^{\sigma}  &  =\left(  -1\right)  ^{\ell\left(
\sigma\right)  }=\left(  -1\right)  ^{\ell\left(  \widehat{\sigma}\right)
}\ \ \ \ \ \ \ \ \ \ \left(  \text{since }\ell\left(  \sigma\right)
=\ell\left(  \widehat{\sigma}\right)  \right) \\
&  =\left(  -1\right)  ^{\widehat{\sigma}}.
\end{align*}
This proves (\ref{pf.lem.laplace.lem.-1}).} and%
\begin{equation}
\prod_{i=1}^{n-1}a_{i,\widehat{\sigma}\left(  i\right)  }=\prod_{i=1}%
^{n-1}a_{i,\sigma\left(  i\right)  } \label{pf.lem.laplace.lem.prod}%
\end{equation}
\footnote{\textit{Proof.} Let $\sigma\in S_{n-1}$. For every $i\in\left\{
1,2,\ldots,n-1\right\}  $, the element $\widehat{\sigma}\left(  i\right)  $ is
well-defined (since $i\in\left\{  1,2,\ldots,n-1\right\}  \subseteq\left\{
1,2,\ldots,n\right\}  $) and satisfies%
\begin{align*}
\widehat{\sigma}\left(  i\right)   &  =%
\begin{cases}
\sigma\left(  i\right)  , & \text{if }i<n;\\
n, & \text{if }i=n
\end{cases}
\ \ \ \ \ \ \ \ \ \ \left(  \text{by the definition of }\widehat{\sigma
}\right) \\
&  =\sigma\left(  i\right)  \ \ \ \ \ \ \ \ \ \ \left(  \text{since }i<n\text{
(since }i\leq n-1\text{ (since }i\in\left\{  1,2,\ldots,n-1\right\}
\text{))}\right)  .
\end{align*}
Hence, for every $i\in\left\{  1,2,\ldots,n-1\right\}  $, we have
$a_{i,\widehat{\sigma}\left(  i\right)  }=a_{i,\sigma\left(  i\right)  }$
(since $\widehat{\sigma}\left(  i\right)  =\sigma\left(  i\right)  $). Thus,
$\prod_{i=1}^{n-1}\underbrace{a_{i,\widehat{\sigma}\left(  i\right)  }%
}_{=a_{i,\sigma\left(  i\right)  }}=\prod_{i=1}^{n-1}a_{i,\sigma\left(
i\right)  }$. This proves (\ref{pf.lem.laplace.lem.prod}).}.
\end{verlong}

\begin{vershort}
Now,
\begin{align*}
&  \underbrace{\sum_{\substack{\sigma\in S_{n};\\\sigma\left(  n\right)  =n}%
}}_{\substack{=\sum_{\sigma\in\left\{  \tau\in S_{n}\ \mid\ \tau\left(
n\right)  =n\right\}  }=\sum_{\sigma\in T}\\\text{(since }\left\{  \tau\in
S_{n}\ \mid\ \tau\left(  n\right)  =n\right\}  =T\text{)}}}\left(  -1\right)
^{\sigma}\prod_{i=1}^{n-1}a_{i,\sigma\left(  i\right)  }\\
&  =\sum_{\sigma\in T}\left(  -1\right)  ^{\sigma}\prod_{i=1}^{n-1}%
a_{i,\sigma\left(  i\right)  }=\sum_{\sigma\in S_{n-1}}\underbrace{\left(
-1\right)  ^{\Phi\left(  \sigma\right)  }\prod_{i=1}^{n-1}a_{i,\left(
\Phi\left(  \sigma\right)  \right)  \left(  i\right)  }}_{\substack{=\left(
-1\right)  ^{\widehat{\sigma}}\prod_{i=1}^{n-1}a_{i,\widehat{\sigma}\left(
i\right)  }\\\text{(since }\Phi\left(  \sigma\right)  =\widehat{\sigma
}\text{)}}}\\
&  \ \ \ \ \ \ \ \ \ \ \left(
\begin{array}
[c]{c}%
\text{here, we have substituted }\Phi\left(  \sigma\right)  \text{ for }%
\sigma\text{ in the sum,}\\
\text{since the map }\Phi:S_{n-1}\rightarrow T\text{ is a bijection}%
\end{array}
\right) \\
&  =\sum_{\sigma\in S_{n-1}}\underbrace{\left(  -1\right)  ^{\widehat{\sigma}%
}}_{\substack{=\left(  -1\right)  ^{\sigma}\\\text{(by
(\ref{pf.lem.laplace.lem.short.-1}))}}}\underbrace{\prod_{i=1}^{n-1}%
a_{i,\widehat{\sigma}\left(  i\right)  }}_{\substack{=\prod_{i=1}%
^{n-1}a_{i,\sigma\left(  i\right)  }\\\text{(by
(\ref{pf.lem.laplace.lem.short.prod}))}}}=\sum_{\sigma\in S_{n-1}}\left(
-1\right)  ^{\sigma}\prod_{i=1}^{n-1}a_{i,\sigma\left(  i\right)  }.
\end{align*}
Compared with%
\begin{align*}
\det\left(  \left(  a_{i,j}\right)  _{1\leq i\leq n-1,\ 1\leq j\leq
n-1}\right)   &  =\sum_{\sigma\in S_{n-1}}\left(  -1\right)  ^{\sigma}%
\prod_{i=1}^{n-1}a_{i,\sigma\left(  i\right)  }\\
&  \ \ \ \ \ \ \ \ \ \ \left(
\begin{array}
[c]{c}%
\text{by (\ref{eq.det.eq.2}), applied to }n-1\text{ and}\\
\left(  a_{i,j}\right)  _{1\leq i\leq n-1,\ 1\leq j\leq n-1}\text{ instead of
}n\text{ and }A
\end{array}
\right)  ,
\end{align*}
this yields%
\[
\sum_{\substack{\sigma\in S_{n};\\\sigma\left(  n\right)  =n}}\left(
-1\right)  ^{\sigma}\prod_{i=1}^{n-1}a_{i,\sigma\left(  i\right)  }%
=\det\left(  \left(  a_{i,j}\right)  _{1\leq i\leq n-1,\ 1\leq j\leq
n-1}\right)  .
\]
This proves Lemma \ref{lem.laplace.lem}. \qedhere

\end{vershort}

\begin{verlong}
Now,%
\begin{align}
&  \underbrace{\sum_{\substack{\sigma\in S_{n};\\\sigma\left(  n\right)  =n}%
}}_{\substack{=\sum_{\sigma\in\left\{  \tau\in S_{n}\ \mid\ \tau\left(
n\right)  =n\right\}  }=\sum_{\sigma\in T}\\\text{(since }\left\{  \tau\in
S_{n}\ \mid\ \tau\left(  n\right)  =n\right\}  =T\text{)}}}\left(  -1\right)
^{\sigma}\prod_{i=1}^{n-1}a_{i,\sigma\left(  i\right)  }\nonumber\\
&  =\sum_{\sigma\in T}\left(  -1\right)  ^{\sigma}\prod_{i=1}^{n-1}%
a_{i,\sigma\left(  i\right)  }=\sum_{\sigma\in S_{n-1}}\underbrace{\left(
-1\right)  ^{\Phi\left(  \sigma\right)  }\prod_{i=1}^{n-1}a_{i,\left(
\Phi\left(  \sigma\right)  \right)  \left(  i\right)  }}_{\substack{=\left(
-1\right)  ^{\widehat{\sigma}}\prod_{i=1}^{n-1}a_{i,\widehat{\sigma}\left(
i\right)  }\\\text{(since }\Phi\left(  \sigma\right)  =\widehat{\sigma
}\text{)}}}\nonumber\\
&  \ \ \ \ \ \ \ \ \ \ \left(
\begin{array}
[c]{c}%
\text{here, we have substituted }\Phi\left(  \sigma\right)  \text{ for }%
\sigma\text{ in the sum,}\\
\text{since the map }\Phi:S_{n-1}\rightarrow T\text{ is a bijection}%
\end{array}
\right) \nonumber\\
&  =\sum_{\sigma\in S_{n-1}}\underbrace{\left(  -1\right)  ^{\widehat{\sigma}%
}}_{\substack{=\left(  -1\right)  ^{\sigma}\\\text{(by
(\ref{pf.lem.laplace.lem.-1}))}}}\underbrace{\prod_{i=1}^{n-1}%
a_{i,\widehat{\sigma}\left(  i\right)  }}_{\substack{=\prod_{i=1}%
^{n-1}a_{i,\sigma\left(  i\right)  }\\\text{(by (\ref{pf.lem.laplace.lem.prod}%
))}}}\nonumber\\
&  =\sum_{\sigma\in S_{n-1}}\left(  -1\right)  ^{\sigma}\prod_{i=1}%
^{n-1}a_{i,\sigma\left(  i\right)  }. \label{pf.lem.laplace.lem.almostthere}%
\end{align}
But (\ref{eq.det.eq.2}) (applied to $n-1$ and $\left(  a_{i,j}\right)  _{1\leq
i\leq n-1,\ 1\leq j\leq n-1}$ instead of $n$ and $A$) yields%
\[
\det\left(  \left(  a_{i,j}\right)  _{1\leq i\leq n-1,\ 1\leq j\leq
n-1}\right)  =\sum_{\sigma\in S_{n-1}}\left(  -1\right)  ^{\sigma}\prod
_{i=1}^{n-1}a_{i,\sigma\left(  i\right)  }.
\]
Compared with (\ref{pf.lem.laplace.lem.almostthere}), this yields%
\[
\sum_{\substack{\sigma\in S_{n};\\\sigma\left(  n\right)  =n}}\left(
-1\right)  ^{\sigma}\prod_{i=1}^{n-1}a_{i,\sigma\left(  i\right)  }%
=\det\left(  \left(  a_{i,j}\right)  _{1\leq i\leq n-1,\ 1\leq j\leq
n-1}\right)  .
\]
This proves Lemma \ref{lem.laplace.lem}.
\end{verlong}
\end{proof}

\begin{proof}
[Proof of Theorem \ref{thm.laplace.pre}.]Every permutation $\sigma\in S_{n}$
satisfying $\sigma\left(  n\right)  \neq n$ satisfies
\begin{equation}
a_{n,\sigma\left(  n\right)  }=0 \label{pf.thm.laplace.pre.1}%
\end{equation}
\footnote{\textit{Proof of (\ref{pf.thm.laplace.pre.1}):} Let $\sigma\in
S_{n}$ be a permutation satisfying $\sigma\left(  n\right)  \neq n$. Since
$\sigma\left(  n\right)  \in\left\{  1,2,\ldots,n\right\}  $ and
$\sigma\left(  n\right)  \neq n$, we have $\sigma\left(  n\right)  \in\left\{
1,2,\ldots,n\right\}  \setminus\left\{  n\right\}  =\left\{  1,2,\ldots
,n-1\right\}  $. Hence, (\ref{eq.thm.laplace.pre.ass}) (applied to
$j=\sigma\left(  n\right)  $) shows that $a_{n,\sigma\left(  n\right)  }=0$,
qed.}.

From (\ref{eq.det.eq.2}), we obtain%
\begin{align*}
\det A  &  =\sum_{\sigma\in S_{n}}\left(  -1\right)  ^{\sigma}%
\underbrace{\prod_{i=1}^{n}a_{i,\sigma\left(  i\right)  }}_{=\left(
\prod_{i=1}^{n-1}a_{i,\sigma\left(  i\right)  }\right)  a_{n,\sigma\left(
n\right)  }}=\sum_{\sigma\in S_{n}}\left(  -1\right)  ^{\sigma}\left(
\prod_{i=1}^{n-1}a_{i,\sigma\left(  i\right)  }\right)  a_{n,\sigma\left(
n\right)  }\\
&  =\sum_{\substack{\sigma\in S_{n};\\\sigma\left(  n\right)  =n}}\left(
-1\right)  ^{\sigma}\left(  \prod_{i=1}^{n-1}a_{i,\sigma\left(  i\right)
}\right)  \underbrace{a_{n,\sigma\left(  n\right)  }}_{\substack{=a_{n,n}%
\\\text{(since }\sigma\left(  n\right)  =n\text{)}}}+\sum_{\substack{\sigma\in
S_{n};\\\sigma\left(  n\right)  \neq n}}\left(  -1\right)  ^{\sigma}\left(
\prod_{i=1}^{n-1}a_{i,\sigma\left(  i\right)  }\right)
\underbrace{a_{n,\sigma\left(  n\right)  }}_{\substack{=0\\\text{(by
(\ref{pf.thm.laplace.pre.1}))}}}\\
&  \ \ \ \ \ \ \ \ \ \ \ \ \ \ \ \ \ \ \ \ \left(
\begin{array}
[c]{c}%
\text{since every }\sigma\in S_{n}\text{ satisfies}\\
\text{either }\sigma\left(  n\right)  =n\text{ or }\sigma\left(  n\right)
\neq n\text{ (but not both)}%
\end{array}
\right) \\
&  =\sum_{\substack{\sigma\in S_{n};\\\sigma\left(  n\right)  =n}}\left(
-1\right)  ^{\sigma}\left(  \prod_{i=1}^{n-1}a_{i,\sigma\left(  i\right)
}\right)  a_{n,n}+\underbrace{\sum_{\substack{\sigma\in S_{n};\\\sigma\left(
n\right)  \neq n}}\left(  -1\right)  ^{\sigma}\left(  \prod_{i=1}%
^{n-1}a_{i,\sigma\left(  i\right)  }\right)  0}_{=0}\\
&  =\sum_{\substack{\sigma\in S_{n};\\\sigma\left(  n\right)  =n}}\left(
-1\right)  ^{\sigma}\left(  \prod_{i=1}^{n-1}a_{i,\sigma\left(  i\right)
}\right)  a_{n,n}=a_{n,n}\cdot\underbrace{\sum_{\substack{\sigma\in
S_{n};\\\sigma\left(  n\right)  =n}}\left(  -1\right)  ^{\sigma}\prod
_{i=1}^{n-1}a_{i,\sigma\left(  i\right)  }}_{\substack{=\det\left(  \left(
a_{i,j}\right)  _{1\leq i\leq n-1,\ 1\leq j\leq n-1}\right)  \\\text{(by Lemma
\ref{lem.laplace.lem})}}}\\
&  =a_{n,n}\cdot\det\left(  \left(  a_{i,j}\right)  _{1\leq i\leq n-1,\ 1\leq
j\leq n-1}\right)  .
\end{align*}
This proves Theorem \ref{thm.laplace.pre}.
\end{proof}

Let us finally state an analogue of Theorem \ref{thm.laplace.pre} in which the
last column (rather than the last row) is required to consist mostly of zeroes:

\begin{corollary}
\label{cor.laplace.pre.col}Let $n$ be a positive integer. Let $A=\left(
a_{i,j}\right)  _{1\leq i\leq n,\ 1\leq j\leq n}$ be an $n\times n$-matrix.
Assume that%
\begin{equation}
a_{i,n}=0\ \ \ \ \ \ \ \ \ \ \text{for every }i\in\left\{  1,2,\ldots
,n-1\right\}  . \label{eq.cor.laplace.pre.col.ass}%
\end{equation}
Then, $\det A=a_{n,n}\cdot\det\left(  \left(  a_{i,j}\right)  _{1\leq i\leq
n-1,\ 1\leq j\leq n-1}\right)  $.
\end{corollary}

\begin{proof}
[Proof of Corollary \ref{cor.laplace.pre.col}.]We have $n-1\in\mathbb{N}$
(since $n$ is a positive integer).

We have $A=\left(  a_{i,j}\right)  _{1\leq i\leq n,\ 1\leq j\leq n}$, and thus
$A^{T}=\left(  a_{j,i}\right)  _{1\leq i\leq n,\ 1\leq j\leq n}$ (by the
definition of the transpose matrix $A^{T}$). Also, for every $j\in\left\{
1,2,\ldots,n-1\right\}  $, we have $a_{j,n}=0$ (by
(\ref{eq.cor.laplace.pre.col.ass}), applied to $i=j$). Thus, Theorem
\ref{thm.laplace.pre} (applied to $A^{T}$ and $a_{j,i}$ instead of $A$ and
$a_{i,j}$) yields
\begin{equation}
\det\left(  A^{T}\right)  =a_{n,n}\cdot\det\left(  \left(  a_{j,i}\right)
_{1\leq i\leq n-1,\ 1\leq j\leq n-1}\right)  .
\label{pf.cor.laplace.pre.col.1}%
\end{equation}

But Exercise \ref{exe.ps4.4} shows that $\det\left(  A^{T}\right)  =\det A$.
Thus, $\det A=\det\left(  A^{T}\right)  $. Also, the definition of the
transpose of a matrix shows that $\left(  \left(  a_{i,j}\right)  _{1\leq
i\leq n-1,\ 1\leq j\leq n-1}\right)  ^{T}=\left(  a_{j,i}\right)  _{1\leq
i\leq n-1,\ 1\leq j\leq n-1}$. Thus,%
\[
\det\left(  \left(  \left(  a_{i,j}\right)  _{1\leq i\leq n-1,\ 1\leq j\leq
n-1}\right)  ^{T}\right)  =\det\left(  \left(  a_{j,i}\right)  _{1\leq i\leq
n-1,\ 1\leq j\leq n-1}\right)  .
\]
Comparing this with%
\[
\det\left(  \left(  \left(  a_{i,j}\right)  _{1\leq i\leq n-1,\ 1\leq j\leq
n-1}\right)  ^{T}\right)  =\det\left(  \left(  a_{i,j}\right)  _{1\leq i\leq
n-1,\ 1\leq j\leq n-1}\right)
\]
(by Exercise \ref{exe.ps4.4}, applied to $n-1$ and $\left(  a_{i,j}\right)
_{1\leq i\leq n-1,\ 1\leq j\leq n-1}$ instead of $n$ and $A$), we obtain%
\[
\det\left(  \left(  a_{j,i}\right)  _{1\leq i\leq n-1,\ 1\leq j\leq
n-1}\right)  =\det\left(  \left(  a_{i,j}\right)  _{1\leq i\leq n-1,\ 1\leq
j\leq n-1}\right)  .
\]

Now,%
\begin{align*}
\det A  &  =\det\left(  A^{T}\right)  =a_{n,n}\cdot\underbrace{\det\left(
\left(  a_{j,i}\right)  _{1\leq i\leq n-1,\ 1\leq j\leq n-1}\right)  }%
_{=\det\left(  \left(  a_{i,j}\right)  _{1\leq i\leq n-1,\ 1\leq j\leq
n-1}\right)  }\ \ \ \ \ \ \ \ \ \ \left(  \text{by
(\ref{pf.cor.laplace.pre.col.1})}\right) \\
&  =a_{n,n}\cdot\det\left(  \left(  a_{i,j}\right)  _{1\leq i\leq n-1,\ 1\leq
j\leq n-1}\right)  .
\end{align*}
This proves Corollary \ref{cor.laplace.pre.col}.
\end{proof}

\subsection{The Vandermonde determinant}

\subsubsection{The statement}

An example for an application of Theorem \ref{thm.laplace.pre} is the famous
\textit{Vandermonde determinant}:

\begin{theorem}
\label{thm.vander-det}Let $n\in\mathbb{N}$. Let $x_{1},x_{2},\ldots,x_{n}$ be
$n$ elements of $\mathbb{K}$. Then:

\textbf{(a)} We have%
\[
\det\left(  \left(  x_{i}^{n-j}\right)  _{1\leq i\leq n,\ 1\leq j\leq
n}\right)  =\prod_{1\leq i<j\leq n}\left(  x_{i}-x_{j}\right)  .
\]

\textbf{(b)} We have%
\[
\det\left(  \left(  x_{j}^{n-i}\right)  _{1\leq i\leq n,\ 1\leq j\leq
n}\right)  =\prod_{1\leq i<j\leq n}\left(  x_{i}-x_{j}\right)  .
\]

\textbf{(c)} We have%
\[
\det\left(  \left(  x_{i}^{j-1}\right)  _{1\leq i\leq n,\ 1\leq j\leq
n}\right)  =\prod_{1\leq j<i\leq n}\left(  x_{i}-x_{j}\right)  .
\]

\textbf{(d)} We have%
\[
\det\left(  \left(  x_{j}^{i-1}\right)  _{1\leq i\leq n,\ 1\leq j\leq
n}\right)  =\prod_{1\leq j<i\leq n}\left(  x_{i}-x_{j}\right)  .
\]

\end{theorem}

\begin{remark}
For $n=4$, the four matrices appearing in Theorem \ref{thm.vander-det} are%
\begin{align*}
\left(  x_{i}^{n-j}\right)  _{1\leq i\leq n,\ 1\leq j\leq n}  &  =\left(
\begin{array}
[c]{cccc}%
x_{1}^{3} & x_{1}^{2} & x_{1} & 1\\
x_{2}^{3} & x_{2}^{2} & x_{2} & 1\\
x_{3}^{3} & x_{3}^{2} & x_{3} & 1\\
x_{4}^{3} & x_{4}^{2} & x_{4} & 1
\end{array}
\right)  ,\\
\left(  x_{j}^{n-i}\right)  _{1\leq i\leq n,\ 1\leq j\leq n}  &  =\left(
\begin{array}
[c]{cccc}%
x_{1}^{3} & x_{2}^{3} & x_{3}^{3} & x_{4}^{3}\\
x_{1}^{2} & x_{2}^{2} & x_{3}^{2} & x_{4}^{2}\\
x_{1} & x_{2} & x_{3} & x_{4}\\
1 & 1 & 1 & 1
\end{array}
\right)  ,\\
\left(  x_{i}^{j-1}\right)  _{1\leq i\leq n,\ 1\leq j\leq n}  &  =\left(
\begin{array}
[c]{cccc}%
1 & x_{1} & x_{1}^{2} & x_{1}^{3}\\
1 & x_{2} & x_{2}^{2} & x_{2}^{3}\\
1 & x_{3} & x_{3}^{2} & x_{3}^{3}\\
1 & x_{4} & x_{4}^{2} & x_{4}^{3}%
\end{array}
\right)  ,\\
\left(  x_{j}^{i-1}\right)  _{1\leq i\leq n,\ 1\leq j\leq n}  &  =\left(
\begin{array}
[c]{cccc}%
1 & 1 & 1 & 1\\
x_{1} & x_{2} & x_{3} & x_{4}\\
x_{1}^{2} & x_{2}^{2} & x_{3}^{2} & x_{4}^{2}\\
x_{1}^{3} & x_{2}^{3} & x_{3}^{3} & x_{4}^{3}%
\end{array}
\right)  .
\end{align*}
It is clear that the second of these four matrices is the transpose of the
first; the fourth is the transpose of the third; and the fourth is obtained
from the second by rearranging the rows in opposite order. Thus, the four
parts of Theorem \ref{thm.vander-det} are rather easily seen to be equivalent.
(We shall prove part \textbf{(a)} and derive the others from it.) Nevertheless
it is useful to have seen them all.
\end{remark}

Theorem \ref{thm.vander-det} is a classical result (known as the
\href{https://en.wikipedia.org/wiki/Vandermonde_matrix}{\textit{Vandermonde
determinant}}, although \href{http://arxiv.org/abs/1204.4716}{it is unclear}
whether it has been proven by Vandermonde): Almost all texts on linear algebra
mention it (or, rather, at least one of its four parts), although some only
prove it in lesser generality. It is a fundamental result that has various
applications to abstract algebra, number theory, coding theory, combinatorics
and numerical mathematics.

Theorem \ref{thm.vander-det} has many known proofs\footnote{For four
combinatorial proofs, see \cite{Gessel-Vand}, \cite[\S 5.3]{Aigner07},
\cite[\S 12.9]{Loehr-BC} and \cite{BenDre-Vand}. (Specifically,
\cite{Gessel-Vand} and \cite{BenDre-Vand} prove Theorem \ref{thm.vander-det}
\textbf{(c)}, whereas \cite[\S 5.3]{Aigner07} and \cite[\S 12.9]{Loehr-BC}
prove Theorem \ref{thm.vander-det} \textbf{(b)}. But as we will see, the four
parts of Theorem \ref{thm.vander-det} are easily seen to be equivalent to each
other.) Gessel's proof from \cite{Gessel-Vand} is also explained in more
detail in \cite{17s-lec8}.}. In this section, I will show two of these proofs.
Another proof (of Theorem \ref{thm.vander-det} \textbf{(c)} only; but as I
said above, the other parts are easily seen to be equivalent) can be found in
Section \ref{sect.sol.det.vdm-pol} or in \cite[Theorem 1]{GriHyp}. Before I
get to the proofs, let me show yet another down-to-earth example of Theorem
\ref{thm.vander-det}:

\begin{example}
\label{exa.vander-det.3}Let $x,y,z\in\mathbb{K}$. Let $A=\left(
\begin{array}
[c]{ccc}%
1 & x & x^{2}\\
1 & y & y^{2}\\
1 & z & z^{2}%
\end{array}
\right)  $. Then, (\ref{eq.det.small.3x3}) shows that%
\begin{align}
\det A  &  =1yz^{2}+xy^{2}\cdot1+x^{2}\cdot1z-1y^{2}z-x\cdot1z^{2}-x^{2}%
y\cdot1\nonumber\\
&  =yz^{2}+xy^{2}+x^{2}z-y^{2}z-xz^{2}-x^{2}y\nonumber\\
&  =yz\left(  z-y\right)  +zx\left(  x-z\right)  +xy\left(  y-x\right)  .
\label{eq.exa.vander-det.3.1}%
\end{align}
On the other hand, Theorem \ref{thm.vander-det} \textbf{(c)} (applied to
$n=3$, $x_{1}=x$, $x_{2}=y$ and $x_{3}=z$) yields $\det A=\left(  y-x\right)
\left(  z-x\right)  \left(  z-y\right)  $. Compared with
(\ref{eq.exa.vander-det.3.1}), this yields%
\begin{equation}
\left(  y-x\right)  \left(  z-x\right)  \left(  z-y\right)  =yz\left(
z-y\right)  +zx\left(  x-z\right)  +xy\left(  y-x\right)  .
\label{eq.exa.vander-det.3.2}%
\end{equation}
You might have encountered this curious identity as a trick of use in contest
problems. When $x,y,z$ are three distinct complex numbers, we can divide
(\ref{eq.exa.vander-det.3.2}) by $\left(  y-x\right)  \left(  z-x\right)
\left(  z-y\right)  $, and obtain%
\[
1=\dfrac{yz}{\left(  y-x\right)  \left(  z-x\right)  }+\dfrac{zx}{\left(
z-y\right)  \left(  x-y\right)  }+\dfrac{xy}{\left(  x-z\right)  \left(
y-z\right)  }.
\]

\end{example}

\subsubsection{A proof by induction}

We now approach the proofs of Theorem \ref{thm.vander-det}. The first proof
has the advantage of demonstrating how Theorem \ref{thm.laplace.pre} can be
used (together with induction) in computing determinants. Before we embark on
this proof, let us see what happens to the determinant of an arbitrary square
matrix if we rearrange the rows in opposite order:

\begin{lemma}
\label{lem.vander-det.lem-rearr}Let $n\in\mathbb{N}$. Let $\left(
a_{i,j}\right)  _{1\leq i\leq n,\ 1\leq j\leq n}$ be an $n\times n$-matrix.
Then,%
\[
\det\left(  \left(  a_{n+1-i,j}\right)  _{1\leq i\leq n,\ 1\leq j\leq
n}\right)  =\left(  -1\right)  ^{n\left(  n-1\right)  /2}\det\left(  \left(
a_{i,j}\right)  _{1\leq i\leq n,\ 1\leq j\leq n}\right)  .
\]

\end{lemma}

\begin{proof}
[Proof of Lemma \ref{lem.vander-det.lem-rearr}.]Let $\left[  n\right]  $
denote the set $\left\{  1,2,\ldots,n\right\}  $. Define a permutation $w_{0}$
in $S_{n}$ as in Exercise \ref{exe.ps4.1c}. In the solution of Exercise
\ref{exe.ps4.1c}, we have shown that $\left(  -1\right)  ^{w_{0}}=\left(
-1\right)  ^{n\left(  n-1\right)  /2}$.

\begin{verlong}
Now, $w_{0}\in S_{n}$. In other words, $w_{0}$ is a permutation of the set
$\left\{  1,2,\ldots,n\right\}  $ (since $S_{n}$ is the set of all
permutations of the set $\left\{  1,2,\ldots,n\right\}  $). In other words,
$w_{0}$ is a permutation of the set $\left[  n\right]  $ (since $\left[
n\right]  =\left\{  1,2,\ldots,n\right\}  $). In other words, $w_{0}$ is a
bijective map $\left[  n\right]  \rightarrow\left[  n\right]  $.
\end{verlong}

Now, we can apply Lemma \ref{lem.det.sigma} \textbf{(a)} to $\left(
a_{i,j}\right)  _{1\leq i\leq n,\ 1\leq j\leq n}$, $w_{0}$ and $\left(
a_{w_{0}\left(  i\right)  ,j}\right)  _{1\leq i\leq n,\ 1\leq j\leq n}$
instead of $B$, $\kappa$ and $B_{\kappa}$. As a result, we obtain%
\begin{align}
\det\left(  \left(  a_{w_{0}\left(  i\right)  ,j}\right)  _{1\leq i\leq
n,\ 1\leq j\leq n}\right)   &  =\underbrace{\left(  -1\right)  ^{w_{0}}%
}_{=\left(  -1\right)  ^{n\left(  n-1\right)  /2}}\cdot\det\left(  \left(
a_{i,j}\right)  _{1\leq i\leq n,\ 1\leq j\leq n}\right) \nonumber\\
&  =\left(  -1\right)  ^{n\left(  n-1\right)  /2}\det\left(  \left(
a_{i,j}\right)  _{1\leq i\leq n,\ 1\leq j\leq n}\right)  .
\label{pf.lem.vander-det.lem-rearr.1}%
\end{align}

\begin{vershort}
But $w_{0}\left(  i\right)  =n+1-i$ for every $i\in\left\{  1,2,\ldots
,n\right\}  $ (by the definition of $w_{0}$). Thus,
(\ref{pf.lem.vander-det.lem-rearr.1}) rewrites as $\det\left(  \left(
a_{n+1-i,j}\right)  _{1\leq i\leq n,\ 1\leq j\leq n}\right)  =\left(
-1\right)  ^{n\left(  n-1\right)  /2}\det\left(  \left(  a_{i,j}\right)
_{1\leq i\leq n,\ 1\leq j\leq n}\right)  $. This proves Lemma
\ref{lem.vander-det.lem-rearr}. \qedhere

\end{vershort}

\begin{verlong}
However, $w_{0}\left(  i\right)  =n+1-i$ for every $i\in\left\{
1,2,\ldots,n\right\}  $ (by the definition of $w_{0}$). Hence, $a_{w_{0}%
\left(  i\right)  ,j}=a_{n+1-i,j}$ for every $i\in\left\{  1,2,\ldots
,n\right\}  $ and $j\in\left\{  1,2,\ldots,n\right\}  $. Thus,%
\[
\det\left(  \left(  \underbrace{a_{w_{0}\left(  i\right)  ,j}}_{=a_{n+1-i,j}%
}\right)  _{1\leq i\leq n,\ 1\leq j\leq n}\right)  =\det\left(  \left(
a_{n+1-i,j}\right)  _{1\leq i\leq n,\ 1\leq j\leq n}\right)  ,
\]
so that%
\begin{align*}
\det\left(  \left(  a_{n+1-i,j}\right)  _{1\leq i\leq n,\ 1\leq j\leq
n}\right)   &  =\det\left(  \left(  a_{w_{0}\left(  i\right)  ,j}\right)
_{1\leq i\leq n,\ 1\leq j\leq n}\right) \\
&  =\left(  -1\right)  ^{n\left(  n-1\right)  /2}\det\left(  \left(
a_{i,j}\right)  _{1\leq i\leq n,\ 1\leq j\leq n}\right)  .
\end{align*}
This proves Lemma \ref{lem.vander-det.lem-rearr}.
\end{verlong}
\end{proof}

\begin{proof}
[First proof of Theorem \ref{thm.vander-det}.]\textbf{(a)} For every
$u\in\left\{  0,1,\ldots,n\right\}  $, let $A_{u}$ be the $u\times u$-matrix
$\left(  x_{i}^{u-j}\right)  _{1\leq i\leq u,\ 1\leq j\leq u}$.

Now, let us show that%
\begin{equation}
\det\left(  A_{u}\right)  =\prod_{1\leq i<j\leq u}\left(  x_{i}-x_{j}\right)
\label{pf.thm.vander-det.a.goal}%
\end{equation}
for every $u\in\left\{  0,1,\ldots,n\right\}  $.

[\textit{Proof of (\ref{pf.thm.vander-det.a.goal}):} We will prove
(\ref{pf.thm.vander-det.a.goal}) by induction over $u$:

\textit{Induction base:} The matrix $A_{0}$ is a $0\times0$-matrix and thus
has determinant $\det\left(  A_{0}\right)  =1$. On the other hand, the product
$\prod_{1\leq i<j\leq0}\left(  x_{i}-x_{j}\right)  $ is an empty product
(i.e., a product of $0$ elements of $\mathbb{K}$) and thus equals $1$ as well.
Hence, both $\det\left(  A_{0}\right)  $ and $\prod_{1\leq i<j\leq0}\left(
x_{i}-x_{j}\right)  $ equal $1$. Thus, $\det\left(  A_{0}\right)
=\prod_{1\leq i<j\leq0}\left(  x_{i}-x_{j}\right)  $. In other words,
(\ref{pf.thm.vander-det.a.goal}) holds for $u=0$. The induction base is thus complete.

\textit{Induction step:} Let $U\in\left\{  1,2,\ldots,n\right\}  $. Assume
that (\ref{pf.thm.vander-det.a.goal}) holds for $u=U-1$. We need to prove that
(\ref{pf.thm.vander-det.a.goal}) holds for $u=U$.

Recall that $A_{U}=\left(  x_{i}^{U-j}\right)  _{1\leq i\leq U,\ 1\leq j\leq
U}$ (by the definition of $A_{U}$).

For every $\left(  i,j\right)  \in\left\{  1,2,\ldots,U\right\}  ^{2}$, define
$b_{i,j}\in\mathbb{K}$ by%
\[
b_{i,j}=%
\begin{cases}
x_{i}^{U-j}-x_{U}x_{i}^{U-j-1}, & \text{if }j<U;\\
1, & \text{if }j=U
\end{cases}
.
\]
Let $B$ be the $U\times U$-matrix $\left(  b_{i,j}\right)  _{1\leq i\leq
U,\ 1\leq j\leq U}$. For example, if $U=4$, then%
\begin{align*}
A_{U}  &  =\left(
\begin{array}
[c]{cccc}%
x_{1}^{3} & x_{1}^{2} & x_{1} & 1\\
x_{2}^{3} & x_{2}^{2} & x_{2} & 1\\
x_{3}^{3} & x_{3}^{2} & x_{3} & 1\\
x_{4}^{3} & x_{4}^{2} & x_{4} & 1
\end{array}
\right)  \ \ \ \ \ \ \ \ \ \ \text{and}\\
B  &  =\left(
\begin{array}
[c]{cccc}%
x_{1}^{3}-x_{4}x_{1}^{2} & x_{1}^{2}-x_{4}x_{1} & x_{1}-x_{4} & 1\\
x_{2}^{3}-x_{4}x_{2}^{2} & x_{2}^{2}-x_{4}x_{2} & x_{2}-x_{4} & 1\\
x_{3}^{3}-x_{4}x_{3}^{2} & x_{3}^{2}-x_{4}x_{3} & x_{3}-x_{4} & 1\\
x_{4}^{3}-x_{4}x_{4}^{2} & x_{4}^{2}-x_{4}x_{4} & x_{4}-x_{4} & 1
\end{array}
\right) \\
&  =\left(
\begin{array}
[c]{cccc}%
x_{1}^{3}-x_{4}x_{1}^{2} & x_{1}^{2}-x_{4}x_{1} & x_{1}-x_{4} & 1\\
x_{2}^{3}-x_{4}x_{2}^{2} & x_{2}^{2}-x_{4}x_{2} & x_{2}-x_{4} & 1\\
x_{3}^{3}-x_{4}x_{3}^{2} & x_{3}^{2}-x_{4}x_{3} & x_{3}-x_{4} & 1\\
0 & 0 & 0 & 1
\end{array}
\right)  .
\end{align*}
We claim that $\det B=\det\left(  A_{U}\right)  $. Indeed, here are two ways
to prove this:

\textit{First proof of }$\det B=\det\left(  A_{U}\right)  $\textit{:} Exercise
\ref{exe.ps4.6k} \textbf{(b)} shows that the determinant of a $U\times
U$-matrix does not change if we subtract a multiple of one of its columns from
another column. Now, let us subtract $x_{U}$ times the $2$-nd column of
$A_{U}$ from the $1$-st column, then subtract $x_{U}$ times the $3$-rd column
of the resulting matrix from the $2$-nd column, and so on, all the way until
we finally subtract $x_{U}$ times the $U$-th column of the matrix from the
$\left(  U-1\right)  $-st column\footnote{So, all in all, we subtract the
$x_{U}$-multiple of each column from its neighbor to its left, but the order
in which we are doing it (namely, from left to right) is important: It means
that the column we are subtracting is unchanged from $A_{U}$. (If we would be
doing these subtractions from right to left instead, then the columns to be
subtracting would be changed by the preceding steps.)}. The resulting matrix
is $B$ (according to our definition of $B$). Thus, $\det B=\det\left(
A_{U}\right)  $ (since our subtractions never change the determinant). This
proves $\det B=\det\left(  A_{U}\right)  $.

\textit{Second proof of }$\det B=\det\left(  A_{U}\right)  $\textit{:} Here is
another way to prove that $\det B=\det\left(  A_{U}\right)  $, with some less handwaving.

For every $\left(  i,j\right)  \in\left\{  1,2,\ldots,U\right\}  ^{2}$, we
define $c_{i,j}\in\mathbb{K}$ by%
\[
c_{i,j}=%
\begin{cases}
1, & \text{if }i=j;\\
-x_{U}, & \text{if }i=j+1;\\
0, & \text{otherwise}%
\end{cases}
.
\]
Let $C$ be the $U\times U$-matrix $\left(  c_{i,j}\right)  _{1\leq i\leq
U,\ 1\leq j\leq U}$.

For example, if $U=4$, then%
\[
C=\left(
\begin{array}
[c]{cccc}%
1 & 0 & 0 & 0\\
-x_{4} & 1 & 0 & 0\\
0 & -x_{4} & 1 & 0\\
0 & 0 & -x_{4} & 1
\end{array}
\right)  .
\]

\begin{vershort}
The matrix $C$ is lower-triangular, and thus Exercise \ref{exe.ps4.3} shows
that its determinant is $\det C=\underbrace{c_{1,1}}_{=1}\underbrace{c_{2,2}%
}_{=1}\cdots\underbrace{c_{U,U}}_{=1}=1$.
\end{vershort}

\begin{verlong}
We have $c_{i,j}=0$ for every $\left(  i,j\right)  \in\left\{  1,2,\ldots
,U\right\}  ^{2}$ satisfying $i<j$\ \ \ \ \footnote{\textit{Proof.} Let
$\left(  i,j\right)  \in\left\{  1,2,\ldots,U\right\}  ^{2}$ be such that
$i<j$. Then, $i\neq j$ (since $i<j$) and $i\neq j+1$ (since $i<j<j+1$). Thus,
we have neither $i=j$ nor $i=j+1$. Now, the definition of $c_{i,j}$ yields
$c_{i,j}=%
\begin{cases}
1, & \text{if }i=j;\\
-x_{U}, & \text{if }i=j+1;\\
0, & \text{otherwise}%
\end{cases}
=0$ (since neither $i=j$ nor $i=j+1$), qed.}. Hence, Exercise \ref{exe.ps4.3}
(applied to $U$, $C$ and $c_{i,j}$ instead of $n$, $A$ and $a_{i,j}$) yields%
\begin{align*}
\det C  &  =c_{1,1}c_{2,2}\cdots c_{U,U}=\prod_{i=1}^{U}\underbrace{c_{i,i}%
}_{\substack{=%
\begin{cases}
1, & \text{if }i=i;\\
-x_{U}, & \text{if }i=i+1;\\
0, & \text{otherwise}%
\end{cases}
\\\text{(by the definition of }c_{i,i}\text{)}}}=\prod_{i=1}^{U}\underbrace{%
\begin{cases}
1, & \text{if }i=i;\\
-x_{U}, & \text{if }i=i+1;\\
0, & \text{otherwise}%
\end{cases}
}_{\substack{=1\\\text{(since }i=i\text{)}}}\\
&  =\prod_{i=1}^{U}1=1.
\end{align*}

\end{verlong}

\begin{vershort}
On the other hand, it is easy to see that $B=A_{U}C$ (check this!). Thus,
Theorem \ref{thm.det(AB)} yields $\det B=\det\left(  A_{U}\right)
\cdot\underbrace{\det C}_{=1}=\det\left(  A_{U}\right)  $. So we have proven
$\det B=\det\left(  A_{U}\right)  $ again.
\end{vershort}

\begin{verlong}
On the other hand, every $\left(  i,j\right)  \in\left\{  1,2,\ldots
,U\right\}  ^{2}$ satisfy%
\begin{equation}
b_{i,j}=\sum_{k=1}^{U}x_{i}^{U-k}c_{k,j} \label{pf.thm.vander-det.B=AUC.pf.1}%
\end{equation}
\footnote{\textit{Proof of (\ref{pf.thm.vander-det.B=AUC.pf.1}):} Let $\left(
i,j\right)  \in\left\{  1,2,\ldots,U\right\}  ^{2}$. We must prove
(\ref{pf.thm.vander-det.B=AUC.pf.1}).
\par
We have $\left(  i,j\right)  \in\left\{  1,2,\ldots,U\right\}  ^{2}$. Thus,
$i\in\left\{  1,2,\ldots,U\right\}  $ and $j\in\left\{  1,2,\ldots,U\right\}
$. We are in one of the following two cases:
\par
\textit{Case 1:} We have $j<U$.
\par
\textit{Case 2:} We have $j\geq U$.
\par
Let us first consider Case 1. In this case, we have $j<U$. Hence,
$j\in\left\{  1,2,\ldots,U-1\right\}  $ (since $j\in\left\{  1,2,\ldots
,U\right\}  $). Now, the definition of $b_{i,j}$ yields%
\begin{equation}
b_{i,j}=%
\begin{cases}
x_{i}^{U-j}-x_{U}x_{i}^{U-j+1}, & \text{if }j<U;\\
1, & \text{if }j=U
\end{cases}
=x_{i}^{U-j}-x_{U}x_{i}^{U-j-1} \label{pf.thm.vander-det.B=AUC.pf.1.pf.1}%
\end{equation}
(since $j<U$). On the other hand, $j\in\left\{  1,2,\ldots,U-1\right\}  $, so
that $j+1\in\left\{  2,3,\ldots,U\right\}  \subseteq\left\{  1,2,\ldots
,U\right\}  $. Thus, $j$ and $j+1$ are two distinct elements of $\left\{
1,2,\ldots,U\right\}  $.
\par
The definition of $c_{j,j}$ yields%
\[
c_{j,j}=%
\begin{cases}
1, & \text{if }j=j;\\
-x_{U}, & \text{if }j=j+1;\\
0, & \text{otherwise}%
\end{cases}
=1\ \ \ \ \ \ \ \ \ \ \left(  \text{since }j=j\right)  .
\]
The definition of $c_{j+1,j}$ yields%
\[
c_{j+1,j}=%
\begin{cases}
1, & \text{if }j+1=j;\\
-x_{U}, & \text{if }j+1=j+1;\\
0, & \text{otherwise}%
\end{cases}
=-x_{U}\ \ \ \ \ \ \ \ \ \ \left(  \text{since }j+1=j+1\right)  .
\]
For every $k\in\left\{  1,2,\ldots,U\right\}  $ satisfying $k\notin\left\{
j,j+1\right\}  $, we have%
\begin{align}
c_{k,j}  &  =%
\begin{cases}
1, & \text{if }k=j;\\
-x_{U}, & \text{if }k=j+1;\\
0, & \text{otherwise}%
\end{cases}
\ \ \ \ \ \ \ \ \ \ \left(  \text{by the definition of }c_{k,j}\right)
\nonumber\\
&  =0\ \ \ \ \ \ \ \ \ \ \left(  \text{since neither }k=j\text{ nor
}k=j+1\text{ (because }k\notin\left\{  j,j+1\right\}  \text{)}\right)
\label{pf.thm.vander-det.B=AUC.pf.1.pf.6}%
\end{align}
\par
Now,%
\begin{align*}
&  \underbrace{\sum_{k=1}^{U}}_{=\sum_{k\in\left\{  1,2,\ldots,U\right\}  }%
}x_{i}^{U-k}c_{k,j}\\
&  =\sum_{k\in\left\{  1,2,\ldots,U\right\}  }x_{i}^{U-k}c_{k,j}\\
&  =x_{i}^{U-j}\underbrace{c_{j,j}}_{=1}+\underbrace{x_{i}^{U-\left(
j+1\right)  }}_{=x_{i}^{U-j-1}}\underbrace{c_{j+1,j}}_{=-x_{U}}+\sum
_{\substack{k\in\left\{  1,2,\ldots,U\right\}  ;\\k\notin\left\{
j,j+1\right\}  }}x_{i}^{U-k}\underbrace{c_{k,j}}_{\substack{=0\\\text{(by
(\ref{pf.thm.vander-det.B=AUC.pf.1.pf.6}))}}}\\
&  \ \ \ \ \ \ \ \ \ \ \left(
\begin{array}
[c]{c}%
\text{here, we have split off the addends for }k=j\text{ and for }k=j+1\text{
from}\\
\text{the sum, because }j\text{ and }j+1\text{ are two distinct elements of
}\left\{  1,2,\ldots,U\right\}
\end{array}
\right) \\
&  =x_{i}^{U-j}1+x_{i}^{U-j-1}\left(  -x_{U}\right)  +\underbrace{\sum
_{\substack{k\in\left\{  1,2,\ldots,U\right\}  ;\\k\notin\left\{
j,j+1\right\}  }}x_{i}^{U-k}0}_{=0}\\
&  =x_{i}^{U-j}1+x_{i}^{U-j-1}\left(  -x_{U}\right)  =x_{i}^{U-j}%
-x_{i}^{U-j-1}x_{U}=x_{i}^{U-j}-x_{U}x_{i}^{U-j-1}.
\end{align*}
Comparing this with (\ref{pf.thm.vander-det.B=AUC.pf.1.pf.1}), we obtain
$b_{i,j}=\sum_{k=1}^{U}x_{i}^{U-k}c_{k,j}$. Thus,
(\ref{pf.thm.vander-det.B=AUC.pf.1}) is proven in Case 1.
\par
Let us now consider Case 2. In this case, we have $j\geq U$. Thus, $j=U$
(since $j\in\left\{  1,2,\ldots,U\right\}  $). Now, the definition of
$b_{i,j}$ yields%
\begin{equation}
b_{i,j}=%
\begin{cases}
x_{i}^{U-j}-x_{U}x_{i}^{U-j+1}, & \text{if }j<U;\\
1, & \text{if }j=U
\end{cases}
=1 \label{pf.thm.vander-det.B=AUC.pf.1.pf.2}%
\end{equation}
(since $j=U$).
\par
Now, let $k\in\left\{  1,2,\ldots,U-1\right\}  $. Then, $k\leq U-1<U=j$ and
thus $k\neq j$. Also, $k<j<j+1$ and thus $k\neq j+1$. Hence, neither $k=j$ nor
$k=j+1$ (since $k\neq j$ and $k\neq j+1$). But the definition of $c_{k,j}$
yields%
\[
c_{k,j}=%
\begin{cases}
1, & \text{if }k=j;\\
-x_{U}, & \text{if }k=j+1;\\
0, & \text{otherwise}%
\end{cases}
=0\ \ \ \ \ \ \ \ \ \ \left(  \text{since neither }k=j\text{ nor
}k=j+1\right)  .
\]
Let us now forget that we fixed $k$. We thus have shown that
\begin{equation}
c_{k,j}=0\ \ \ \ \ \ \ \ \ \ \text{for every }k\in\left\{  1,2,\ldots
,U-1\right\}  . \label{pf.thm.vander-det.B=AUC.pf.1.pf.3}%
\end{equation}
Now,
\begin{align*}
&  \sum_{k=1}^{U}x_{i}^{U-k}c_{k,j}=\sum_{k=1}^{U-1}x_{i}^{U-k}%
\underbrace{c_{k,j}}_{\substack{=0\\\text{(by
(\ref{pf.thm.vander-det.B=AUC.pf.1.pf.3}))}}}+\underbrace{x_{i}^{U-U}}%
_{=x_{i}^{0}=1}\underbrace{c_{U,j}}_{\substack{=%
\begin{cases}
1, & \text{if }U=j;\\
-x_{U}, & \text{if }U=j+1;\\
0, & \text{otherwise}%
\end{cases}
\\\text{(by the definition of }c_{U,j}\text{)}}}\\
&  =\underbrace{\sum_{k=1}^{U-1}x_{i}^{U-k}0}_{=0}+1\underbrace{%
\begin{cases}
1, & \text{if }U=j;\\
-x_{U}, & \text{if }U=j+1;\\
0, & \text{otherwise}%
\end{cases}
}_{\substack{=1\\\text{(since }U=j\text{)}}}=0+1\cdot1=1.
\end{align*}
Compared with (\ref{pf.thm.vander-det.B=AUC.pf.1.pf.2}), this yields
$b_{i,j}=\sum_{k=1}^{U}x_{i}^{U-k}c_{k,j}$. Thus,
(\ref{pf.thm.vander-det.B=AUC.pf.1}) is proven in Case 2.
\par
We have now proven (\ref{pf.thm.vander-det.B=AUC.pf.1}) in each of the two
Cases 1 and 2. Since these two Cases cover all possibilities, this yields that
(\ref{pf.thm.vander-det.B=AUC.pf.1}) always holds. Qed.}. Now,%
\[
B=\left(  \underbrace{b_{i,j}}_{\substack{=\sum_{k=1}^{U}x_{i}^{U-k}%
c_{k,j}\\\text{(by (\ref{pf.thm.vander-det.B=AUC.pf.1}))}}}\right)  _{1\leq
i\leq U,\ 1\leq j\leq U}=\left(  \sum_{k=1}^{U}x_{i}^{U-k}c_{k,j}\right)
_{1\leq i\leq U,\ 1\leq j\leq U}.
\]
Compared with%
\begin{align*}
A_{U}C  &  =\left(  \sum_{k=1}^{U}x_{i}^{U-k}c_{k,j}\right)  _{1\leq i\leq
U,\ 1\leq j\leq U}\\
&  \ \ \ \ \ \ \ \ \ \ \left(
\begin{array}
[c]{c}%
\text{by the definition of the product }A_{U}C\text{,}\\
\text{since }A_{U}=\left(  x_{i}^{U-j}\right)  _{1\leq i\leq U,\ 1\leq j\leq
U}\text{ and }C=\left(  c_{i,j}\right)  _{1\leq i\leq U,\ 1\leq j\leq U}%
\end{array}
\right)  ,
\end{align*}
this yields $B=A_{U}C$. Hence,%
\begin{align*}
\det\underbrace{B}_{=A_{U}C}  &  =\det\left(  A_{U}C\right)  =\det\left(
A_{U}\right)  \cdot\underbrace{\det C}_{=1}\\
&  \ \ \ \ \ \ \ \ \ \ \left(
\begin{array}
[c]{c}%
\text{by Theorem \ref{thm.det(AB)}, applied to }U\text{, }A_{U}\text{ and }C\\
\text{instead of }n\text{, }A\text{ and }B
\end{array}
\right) \\
&  =\det\left(  A_{U}\right)  .
\end{align*}
Thus, $\det B=\det\left(  A_{U}\right)  $ is proven again.
\end{verlong}

[\textit{Remark:} It is instructive to compare the two proofs of $\det
B=\det\left(  A_{U}\right)  $ given above. They are close kin, although they
might look different at first. In the first proof, we argued that $B$ can be
obtained from $A_{U}$ by subtracting multiples of some columns from others; in
the second, we argued that $B=A_{U}C$ for a specific lower-triangular matrix
$C$. But a look at the matrix $C$ makes it clear that multiplying a $U\times
U$-matrix with $C$ on the right (i.e., transforming a $U\times U$-matrix $X$
into the matrix $XC$) is tantamount to subtracting multiples of some columns
from others, in the way we did it to $A_{U}$ to obtain $B$. So the main
difference between the two proofs is that the first proof used a step-by-step
procedure to obtain $B$ from $A_{U}$, whereas the second proof obtained $B$
from $A_{U}$ by a single-step operation (namely, multiplication by a matrix
$C$).]

Next, we observe that for every $j\in\left\{  1,2,\ldots,U-1\right\}  $, we
have%
\begin{align*}
b_{U,j}  &  =%
\begin{cases}
x_{U}^{U-j}-x_{U}x_{U}^{U-j-1}, & \text{if }j<U;\\
1, & \text{if }j=U
\end{cases}
\ \ \ \ \ \ \ \ \ \ \left(  \text{by the definition of }b_{U,j}\right) \\
&  =x_{U}^{U-j}-\underbrace{x_{U}x_{U}^{U-j-1}}_{=x_{U}^{\left(  U-j-1\right)
+1}=x_{U}^{U-j}}\ \ \ \ \ \ \ \ \ \ \left(  \text{since }j<U\text{ (since
}j\in\left\{  1,2,\ldots,U-1\right\}  \text{)}\right) \\
&  =x_{U}^{U-j}-x_{U}^{U-j}=0.
\end{align*}
Hence, Theorem \ref{thm.laplace.pre} (applied to $U$, $B$ and $b_{i,j}$
instead of $n$, $A$ and $a_{i,j}$) yields%
\begin{equation}
\det B=b_{U,U}\cdot\det\left(  \left(  b_{i,j}\right)  _{1\leq i\leq
U-1,\ 1\leq j\leq U-1}\right)  . \label{pf.thm.vander-det.detB=prod}%
\end{equation}
Let $B^{\prime}$ denote the $\left(  U-1\right)  \times\left(  U-1\right)
$-matrix $\left(  b_{i,j}\right)  _{1\leq i\leq U-1,\ 1\leq j\leq U-1}$.

The definition of $b_{U,U}$ yields%
\begin{align*}
b_{U,U}  &  =%
\begin{cases}
x_{U}^{U-U}-x_{U}x_{U}^{U-U-1}, & \text{if }U<U;\\
1, & \text{if }U=U
\end{cases}
\ \ \ \ \ \ \ \ \ \ \left(  \text{by the definition of }b_{U,U}\right) \\
&  =1\ \ \ \ \ \ \ \ \ \ \left(  \text{since }U=U\right)  .
\end{align*}
Thus, (\ref{pf.thm.vander-det.detB=prod}) becomes%
\[
\det B=\underbrace{b_{U,U}}_{=1}\cdot\det\left(  \underbrace{\left(
b_{i,j}\right)  _{1\leq i\leq U-1,\ 1\leq j\leq U-1}}_{=B^{\prime}}\right)
=\det\left(  B^{\prime}\right)  .
\]
Compared with $\det B=\det\left(  A_{U}\right)  $, this yields%
\begin{equation}
\det\left(  A_{U}\right)  =\det\left(  B^{\prime}\right)  .
\label{pf.thm.vander-det.detAU=detB'}%
\end{equation}

Now, let us take a closer look at $B^{\prime}$. Indeed, every $\left(
i,j\right)  \in\left\{  1,2,\ldots,U-1\right\}  ^{2}$ satisfies%
\begin{align}
b_{i,j}  &  =%
\begin{cases}
x_{i}^{U-j}-x_{U}x_{i}^{U-j-1}, & \text{if }j<U;\\
1, & \text{if }j=U
\end{cases}
\ \ \ \ \ \ \ \ \ \ \left(  \text{by the definition of }b_{i,j}\right)
\nonumber\\
&  =\underbrace{x_{i}^{U-j}}_{=x_{i}x_{i}^{U-j-1}}-x_{U}x_{i}^{U-j-1}%
\ \ \ \ \ \ \ \ \ \ \left(
\begin{array}
[c]{c}%
\text{since }j<U\text{ (since }j\in\left\{  1,2,\ldots,U-1\right\} \\
\text{(since }\left(  i,j\right)  \in\left\{  1,2,\ldots,U-1\right\}
^{2}\text{))}%
\end{array}
\right) \nonumber\\
&  =x_{i}x_{i}^{U-j-1}-x_{U}x_{i}^{U-j-1}=\left(  x_{i}-x_{U}\right)
\underbrace{x_{i}^{U-j-1}}_{=x_{i}^{\left(  U-1\right)  -j}}=\left(
x_{i}-x_{U}\right)  x_{i}^{\left(  U-1\right)  -j}.
\label{pf.thm.vander-det.bij-small}%
\end{align}
Hence,%
\begin{equation}
B^{\prime}=\left(  \underbrace{b_{i,j}}_{\substack{=\left(  x_{i}%
-x_{U}\right)  x_{i}^{\left(  U-1\right)  -j}\\\text{(by
(\ref{pf.thm.vander-det.bij-small}))}}}\right)  _{1\leq i\leq U-1,\ 1\leq
j\leq U-1}=\left(  \left(  x_{i}-x_{U}\right)  x_{i}^{\left(  U-1\right)
-j}\right)  _{1\leq i\leq U-1,\ 1\leq j\leq U-1}. \label{pf.thm.vander-det.B'}%
\end{equation}
On the other hand, the definition of $A_{U-1}$ yields
\begin{equation}
A_{U-1}=\left(  x_{i}^{\left(  U-1\right)  -j}\right)  _{1\leq i\leq
U-1,\ 1\leq j\leq U-1}. \label{pf.thm.vander-det.AU-1}%
\end{equation}

Now, we claim that%
\begin{equation}
\det\left(  B^{\prime}\right)  =\det\left(  A_{U-1}\right)  \cdot\prod
_{i=1}^{U-1}\left(  x_{i}-x_{U}\right)  . \label{pf.thm.vander-det.detB'=}%
\end{equation}
Indeed, here are two ways to prove this:

\textit{First proof of (\ref{pf.thm.vander-det.detB'=}):} Comparing the
formulas (\ref{pf.thm.vander-det.B'}) and (\ref{pf.thm.vander-det.AU-1}), we
see that the matrix $B^{\prime}$ is obtained from the matrix $A_{U-1}$ by
multiplying the first row by $x_{1}-x_{U}$, the second row by $x_{2}-x_{U}$,
and so on, and finally the $\left(  U-1\right)  $-st row by $x_{U-1}-x_{U}$.
But every time we multiply a row of a $\left(  U-1\right)  \times\left(
U-1\right)  $-matrix by some scalar $\lambda\in\mathbb{K}$, the determinant of
the matrix gets multiplied by $\lambda$ (because of Exercise \ref{exe.ps4.6}
\textbf{(g)}). Hence, the determinant of $B^{\prime}$ is obtained from that of
$A_{U-1}$ by first multiplying by $x_{1}-x_{U}$, then multiplying by
$x_{2}-x_{U}$, and so on, and finally multiplying with $x_{U-1}-x_{U}$. In
other words,%
\[
\det\left(  B^{\prime}\right)  =\det\left(  A_{U-1}\right)  \cdot\prod
_{i=1}^{U-1}\left(  x_{i}-x_{U}\right)  .
\]
This proves (\ref{pf.thm.vander-det.detB'=}).

\textit{Second proof of (\ref{pf.thm.vander-det.detB'=}):} For every $\left(
i,j\right)  \in\left\{  1,2,\ldots,U-1\right\}  ^{2}$, we define $d_{i,j}%
\in\mathbb{K}$ by%
\[
d_{i,j}=%
\begin{cases}
x_{i}-x_{U}, & \text{if }i=j;\\
0, & \text{otherwise}%
\end{cases}
.
\]
Let $D$ be the $\left(  U-1\right)  \times\left(  U-1\right)  $-matrix
$\left(  d_{i,j}\right)  _{1\leq i\leq U-1,\ 1\leq j\leq U-1}$.

For example, if $U=4$, then%
\[
D=\left(
\begin{array}
[c]{ccc}%
x_{1}-x_{4} & 0 & 0\\
0 & x_{2}-x_{4} & 0\\
0 & 0 & x_{3}-x_{4}%
\end{array}
\right)  .
\]

\begin{vershort}
The matrix $D$ is lower-triangular (actually, diagonal\footnote{A square
matrix $E=\left(  e_{i,j}\right)  _{1\leq i\leq n,\ 1\leq j\leq n}$ is said to
be \textit{diagonal} if every $\left(  i,j\right)  \in\left\{  1,2,\ldots
,n\right\}  ^{2}$ satisfying $i\neq j$ satisfies $e_{i,j}=0$. In other words,
a square matrix is said to be \textit{diagonal} if it is both upper-triangular
and lower-triangular.}), and thus Exercise \ref{exe.ps4.3} shows that its
determinant is%
\[
\det D=\left(  x_{1}-x_{U}\right)  \left(  x_{2}-x_{U}\right)  \cdots\left(
x_{U-1}-x_{U}\right)  =\prod_{i=1}^{U-1}\left(  x_{i}-x_{U}\right)  .
\]

\end{vershort}

\begin{verlong}
We have $d_{i,j}=0$ for every $\left(  i,j\right)  \in\left\{  1,2,\ldots
,U-1\right\}  ^{2}$ satisfying $i<j$\ \ \ \ \footnote{\textit{Proof.} Let
$\left(  i,j\right)  \in\left\{  1,2,\ldots,U-1\right\}  ^{2}$ be such that
$i<j$. Then, we don't have $i=j$ (since $i<j$). Now, the definition of
$d_{i,j}$ yields $d_{i,j}=%
\begin{cases}
x_{i}-x_{U}, & \text{if }i=j;\\
0, & \text{otherwise}%
\end{cases}
=0$ (since we don't have $i=j$), qed.}. Hence, Exercise \ref{exe.ps4.3}
(applied to $U-1$, $D$ and $d_{i,j}$ instead of $n$, $A$ and $a_{i,j}$) yields%
\begin{align*}
\det D  &  =d_{1,1}d_{2,2}\cdots d_{U-1,U-1}=\prod_{i=1}^{U-1}%
\underbrace{d_{i,i}}_{\substack{=%
\begin{cases}
x_{i}-x_{U}, & \text{if }i=i;\\
0, & \text{otherwise}%
\end{cases}
\\\text{(by the definition of }d_{i,i}\text{)}}}\\
&  =\prod_{i=1}^{U-1}\underbrace{%
\begin{cases}
x_{i}-x_{U}, & \text{if }i=i;\\
0, & \text{otherwise}%
\end{cases}
}_{\substack{=x_{i}-x_{U}\\\text{(since }i=i\text{)}}}=\prod_{i=1}%
^{U-1}\left(  x_{i}-x_{U}\right)  .
\end{align*}

\end{verlong}

\begin{vershort}
On the other hand, it is easy to see that $B^{\prime}=DA_{U-1}$ (check this!).
Thus, Theorem \ref{thm.det(AB)} yields
\[
\det\left(  B^{\prime}\right)  =\det D\cdot\det\left(  A_{U-1}\right)
=\det\left(  A_{U-1}\right)  \cdot\underbrace{\det D}_{=\prod_{i=1}%
^{U-1}\left(  x_{i}-x_{U}\right)  }=\det\left(  A_{U-1}\right)  \cdot
\prod_{i=1}^{U-1}\left(  x_{i}-x_{U}\right)  .
\]
Thus, (\ref{pf.thm.vander-det.detB'=}) is proven again.
\end{vershort}

\begin{verlong}
On the other hand, every $\left(  i,j\right)  \in\left\{  1,2,\ldots
,U-1\right\}  ^{2}$ satisfy%
\begin{equation}
b_{i,j}=\sum_{k=1}^{U-1}d_{i,k}x_{k}^{\left(  U-1\right)  -j}
\label{pf.thm.vander-det.B'=DAU-1.pf.1}%
\end{equation}
\footnote{\textit{Proof of (\ref{pf.thm.vander-det.B'=DAU-1.pf.1}):} Let
$\left(  i,j\right)  \in\left\{  1,2,\ldots,U-1\right\}  ^{2}$. We must prove
(\ref{pf.thm.vander-det.B'=DAU-1.pf.1}).
\par
We have $\left(  i,j\right)  \in\left\{  1,2,\ldots,U-1\right\}  ^{2}$. Thus,
$i\in\left\{  1,2,\ldots,U-1\right\}  $ and $j\in\left\{  1,2,\ldots
,U-1\right\}  $. Now,%
\begin{align*}
&  \underbrace{\sum_{k=1}^{U-1}}_{=\sum_{k\in\left\{  1,2,\ldots,U-1\right\}
}}\underbrace{d_{i,k}}_{\substack{=%
\begin{cases}
x_{i}-x_{U}, & \text{if }i=k;\\
0, & \text{otherwise}%
\end{cases}
\\\text{(by the definition of }d_{i,k}\text{)}}}x_{k}^{\left(  U-1\right)
-j}\\
&  =\sum_{k\in\left\{  1,2,\ldots,U-1\right\}  }%
\begin{cases}
x_{i}-x_{U}, & \text{if }i=k;\\
0, & \text{otherwise}%
\end{cases}
x_{k}^{\left(  U-1\right)  -j}\\
&  =\underbrace{%
\begin{cases}
x_{i}-x_{U}, & \text{if }i=i;\\
0, & \text{otherwise}%
\end{cases}
}_{\substack{=x_{i}-x_{U}\\\text{(since }i=i\text{)}}}x_{i}^{\left(
U-1\right)  -j}+\sum_{\substack{k\in\left\{  1,2,\ldots,U-1\right\}  ;\\k\neq
i}}\underbrace{%
\begin{cases}
x_{i}-x_{U}, & \text{if }i=k;\\
0, & \text{otherwise}%
\end{cases}
}_{\substack{=0\\\text{(since we don't have }i=k\\\text{(since }i\neq k\text{
(since }k\neq i\text{)))}}}x_{k}^{\left(  U-1\right)  -j}\\
&  \ \ \ \ \ \ \ \ \ \ \left(  \text{here, we have split off the addend for
}k=i\text{ from the sum}\right) \\
&  =\left(  x_{i}-x_{U}\right)  x_{i}^{\left(  U-1\right)  -j}%
+\underbrace{\sum_{\substack{k\in\left\{  1,2,\ldots,U-1\right\}  ;\\k\neq
i}}0x_{k}^{\left(  U-1\right)  -j}}_{=0}=\left(  x_{i}-x_{U}\right)
x_{i}^{\left(  U-1\right)  -j}=b_{i,j}\ \ \ \ \ \ \ \ \ \ \left(  \text{by
(\ref{pf.thm.vander-det.bij-small})}\right)  ,
\end{align*}
and therefore (\ref{pf.thm.vander-det.B'=DAU-1.pf.1}) is proven.}. Now,%
\[
B^{\prime}=\left(  \underbrace{b_{i,j}}_{\substack{=\sum_{k=1}^{U-1}%
d_{i,k}x_{k}^{\left(  U-1\right)  -j}\\\text{(by
(\ref{pf.thm.vander-det.B'=DAU-1.pf.1}))}}}\right)  _{1\leq i\leq U-1,\ 1\leq
j\leq U-1}=\left(  \sum_{k=1}^{U-1}d_{i,k}x_{k}^{\left(  U-1\right)
-j}\right)  _{1\leq i\leq U-1,\ 1\leq j\leq U-1}.
\]
Compared with%
\begin{align*}
DA_{U-1}  &  =\left(  \sum_{k=1}^{U-1}d_{i,k}x_{k}^{\left(  U-1\right)
-j}\right)  _{1\leq i\leq U-1,\ 1\leq j\leq U-1}\\
&  \ \ \ \ \ \ \ \ \ \ \left(
\begin{array}
[c]{c}%
\text{by the definition of the product }DA_{U-1}\text{,}\\
\text{since }D=\left(  d_{i,j}\right)  _{1\leq i\leq U-1,\ 1\leq j\leq U-1}\\
\text{and }A_{U-1}=\left(  x_{i}^{\left(  U-1\right)  -j}\right)  _{1\leq
i\leq U-1,\ 1\leq j\leq U-1}%
\end{array}
\right)  ,
\end{align*}
this yields $B^{\prime}=DA_{U-1}$. Hence,%
\begin{align*}
\det\left(  \underbrace{B^{\prime}}_{=DA_{U-1}}\right)   &  =\det\left(
DA_{U-1}\right)  =\det D\cdot\det\left(  A_{U-1}\right) \\
&  \ \ \ \ \ \ \ \ \ \ \left(
\begin{array}
[c]{c}%
\text{by Theorem \ref{thm.det(AB)}, applied to }U-1\text{, }D\text{ and
}A_{U-1}\\
\text{instead of }n\text{, }A\text{ and }B
\end{array}
\right) \\
&  =\det\left(  A_{U-1}\right)  \cdot\underbrace{\det D}_{=\prod_{i=1}%
^{U-1}\left(  x_{i}-x_{U}\right)  }=\det\left(  A_{U-1}\right)  \cdot
\prod_{i=1}^{U-1}\left(  x_{i}-x_{U}\right)  .
\end{align*}
Thus, (\ref{pf.thm.vander-det.detB'=}) is proven again.
\end{verlong}

[\textit{Remark:} Again, our two proofs of (\ref{pf.thm.vander-det.detB'=})
are closely related: the first one reveals $B^{\prime}$ as the result of a
step-by-step process applied to $A_{U-1}$, while the second shows how
$B^{\prime}$ can be obtained from $A_{U-1}$ by a single multiplication.
However, here (in contrast to the proofs of $\det B=\det\left(  A_{U}\right)
$), the step-by-step process involves transforming rows (not columns), and the
multiplication is a multiplication from the left (we have $B^{\prime}%
=DA_{U-1}$, not $B^{\prime}=A_{U-1}D$).]

Now, (\ref{pf.thm.vander-det.detAU=detB'}) becomes%
\begin{equation}
\det\left(  A_{U}\right)  =\det\left(  B^{\prime}\right)  =\det\left(
A_{U-1}\right)  \cdot\prod_{i=1}^{U-1}\left(  x_{i}-x_{U}\right)  .
\label{pf.thm.vander-det.detAU=}%
\end{equation}

But we have assumed that (\ref{pf.thm.vander-det.a.goal}) holds for $u=U-1$.
In other words,%
\[
\det\left(  A_{U-1}\right)  =\underbrace{\prod_{1\leq i<j\leq U-1}}%
_{=\prod_{j=1}^{U-1}\ \ \prod_{i=1}^{j-1}}\left(  x_{i}-x_{j}\right)
=\prod_{j=1}^{U-1}\ \ \prod_{i=1}^{j-1}\left(  x_{i}-x_{j}\right)  .
\]
Hence, (\ref{pf.thm.vander-det.detAU=}) yields%
\begin{align*}
\det\left(  A_{U}\right)   &  =\underbrace{\det\left(  A_{U-1}\right)
}_{=\prod_{j=1}^{U-1}\ \ \prod_{i=1}^{j-1}\left(  x_{i}-x_{j}\right)  }%
\cdot\prod_{i=1}^{U-1}\left(  x_{i}-x_{U}\right) \\
&  =\left(  \prod_{j=1}^{U-1}\ \ \prod_{i=1}^{j-1}\left(  x_{i}-x_{j}\right)
\right)  \cdot\prod_{i=1}^{U-1}\left(  x_{i}-x_{U}\right)  .
\end{align*}
Compared with%
\begin{align*}
\underbrace{\prod_{1\leq i<j\leq U}}_{=\prod_{j=1}^{U}\ \ \prod_{i=1}^{j-1}%
}\left(  x_{i}-x_{j}\right)   &  =\prod_{j=1}^{U}\ \ \prod_{i=1}^{j-1}\left(
x_{i}-x_{j}\right)  =\left(  \prod_{j=1}^{U-1}\ \ \prod_{i=1}^{j-1}\left(
x_{i}-x_{j}\right)  \right)  \cdot\prod_{i=1}^{U-1}\left(  x_{i}-x_{U}\right)
\\
&  \ \ \ \ \ \ \ \ \ \ \left(  \text{here, we have split off the factor for
}j=U\text{ from the product}\right)  ,
\end{align*}
this yields $\det\left(  A_{U}\right)  =\prod_{1\leq i<j\leq U}\left(
x_{i}-x_{j}\right)  $. In other words, (\ref{pf.thm.vander-det.a.goal}) holds
for $u=U$. This completes the induction step.

Now, (\ref{pf.thm.vander-det.a.goal}) is proven by induction.]

Hence, we can apply (\ref{pf.thm.vander-det.a.goal}) to $u=n$. As the result,
we obtain $\det\left(  A_{n}\right)  =\prod_{1\leq i<j\leq n}\left(
x_{i}-x_{j}\right)  $. Since $A_{n}=\left(  x_{i}^{n-j}\right)  _{1\leq i\leq
n,\ 1\leq j\leq n}$ (by the definition of $A_{n}$), this rewrites as
$\det\left(  \left(  x_{i}^{n-j}\right)  _{1\leq i\leq n,\ 1\leq j\leq
n}\right)  =\prod_{1\leq i<j\leq n}\left(  x_{i}-x_{j}\right)  $. This proves
Theorem \ref{thm.vander-det} \textbf{(a)}.

\textbf{(b)} The definition of the transpose of a matrix yields $\left(
\left(  x_{j}^{n-i}\right)  _{1\leq i\leq n,\ 1\leq j\leq n}\right)
^{T}=\left(  x_{i}^{n-j}\right)  _{1\leq i\leq n,\ 1\leq j\leq n}$. Hence,%
\[
\det\left(  \underbrace{\left(  \left(  x_{j}^{n-i}\right)  _{1\leq i\leq
n,\ 1\leq j\leq n}\right)  ^{T}}_{=\left(  x_{i}^{n-j}\right)  _{1\leq i\leq
n,\ 1\leq j\leq n}}\right)  =\det\left(  \left(  x_{i}^{n-j}\right)  _{1\leq
i\leq n,\ 1\leq j\leq n}\right)  =\prod_{1\leq i<j\leq n}\left(  x_{i}%
-x_{j}\right)
\]
(by Theorem \ref{thm.vander-det} \textbf{(a)}). Compared with%
\[
\det\left(  \left(  \left(  x_{j}^{n-i}\right)  _{1\leq i\leq n,\ 1\leq j\leq
n}\right)  ^{T}\right)  =\det\left(  \left(  x_{j}^{n-i}\right)  _{1\leq i\leq
n,\ 1\leq j\leq n}\right)
\]
(by Exercise \ref{exe.ps4.4}, applied to $A=\left(  x_{j}^{n-i}\right)
_{1\leq i\leq n,\ 1\leq j\leq n}$), this yields%
\[
\det\left(  \left(  x_{j}^{n-i}\right)  _{1\leq i\leq n,\ 1\leq j\leq
n}\right)  =\prod_{1\leq i<j\leq n}\left(  x_{i}-x_{j}\right)  .
\]
This proves Theorem \ref{thm.vander-det} \textbf{(b)}.

\textbf{(d)} Applying Lemma \ref{lem.vander-det.lem-rearr} to $a_{i,j}%
=x_{j}^{n-i}$, we obtain%
\begin{align*}
\det\left(  \left(  x_{j}^{n-\left(  n+1-i\right)  }\right)  _{1\leq i\leq
n,\ 1\leq j\leq n}\right)   &  =\left(  -1\right)  ^{n\left(  n-1\right)
/2}\underbrace{\det\left(  \left(  x_{j}^{n-i}\right)  _{1\leq i\leq n,\ 1\leq
j\leq n}\right)  }_{\substack{=\prod_{1\leq i<j\leq n}\left(  x_{i}%
-x_{j}\right)  \\\text{(by Theorem \ref{thm.vander-det} \textbf{(b)})}}}\\
&  =\left(  -1\right)  ^{n\left(  n-1\right)  /2}\prod_{1\leq i<j\leq
n}\left(  x_{i}-x_{j}\right)  .
\end{align*}
This rewrites as%
\begin{equation}
\det\left(  \left(  x_{j}^{i-1}\right)  _{1\leq i\leq n,\ 1\leq j\leq
n}\right)  =\left(  -1\right)  ^{n\left(  n-1\right)  /2}\prod_{1\leq i<j\leq
n}\left(  x_{i}-x_{j}\right)  \label{pf.thm.vander-det.d.0}%
\end{equation}
(since every $\left(  i,j\right)  \in\left\{  1,2,\ldots,n\right\}  ^{2}$
satisfies $x_{j}^{n-\left(  n+1-i\right)  }=x_{j}^{i-1}$).

\begin{vershort}
Now, in the solution to Exercise \ref{exe.ps4.1c}, we have shown that the
number of all pairs $\left(  i,j\right)  $ of integers satisfying $1\leq
i<j\leq n$ is $n\left(  n-1\right)  /2$. In other words,%
\begin{equation}
\left(  \text{the number of all }\left(  i,j\right)  \in\left\{
1,2,\ldots,n\right\}  ^{2}\text{ such that }i<j\right)  =n\left(  n-1\right)
/2. \label{pf.thm.vander-det.d.short.1}%
\end{equation}
Now,
\begin{align*}
\prod_{1\leq j<i\leq n}\left(  x_{i}-x_{j}\right)   &  =\prod_{1\leq i<j\leq
n}\underbrace{\left(  x_{j}-x_{i}\right)  }_{=\left(  -1\right)  \left(
x_{i}-x_{j}\right)  }\ \ \ \ \ \ \ \ \ \ \left(
\begin{array}
[c]{c}%
\text{here, we renamed the index }\left(  j,i\right) \\
\text{as }\left(  i,j\right)  \text{ in the product}%
\end{array}
\right) \\
&  =\prod_{1\leq i<j\leq n}\left(  \left(  -1\right)  \left(  x_{i}%
-x_{j}\right)  \right) \\
&  =\underbrace{\left(  -1\right)  ^{\left(  \text{the number of all }\left(
i,j\right)  \in\left\{  1,2,\ldots,n\right\}  ^{2}\text{ such that
}i<j\right)  }}_{\substack{=\left(  -1\right)  ^{n\left(  n-1\right)
/2}\\\text{(by (\ref{pf.thm.vander-det.d.short.1}))}}}\prod_{1\leq i<j\leq
n}\left(  x_{i}-x_{j}\right) \\
&  =\left(  -1\right)  ^{n\left(  n-1\right)  /2}\prod_{1\leq i<j\leq
n}\left(  x_{i}-x_{j}\right)  .
\end{align*}
Compared with (\ref{pf.thm.vander-det.d.0}), this yields $\det\left(  \left(
x_{j}^{i-1}\right)  _{1\leq i\leq n,\ 1\leq j\leq n}\right)  =\prod_{1\leq
j<i\leq n}\left(  x_{i}-x_{j}\right)  $. This proves Theorem
\ref{thm.vander-det} \textbf{(d)}.
\end{vershort}

\begin{verlong}
Let $G=\left\{  \left(  i,j\right)  \in\mathbb{Z}^{2}\ \mid\ 1\leq i<j\leq
n\right\}  $. Then, $\left\vert G\right\vert =n\left(  n-1\right)  /2$. (This
was proven in the solution to Exercise \ref{exe.ps4.1c}.) Now,%
\begin{align}
\prod_{1\leq j<i\leq n}\left(  x_{i}-x_{j}\right)   &  =\underbrace{\prod
_{1\leq i<j\leq n}}_{\substack{=\prod_{\substack{\left(  i,j\right)
\in\mathbb{Z}^{2};\\1\leq i<j\leq n}}=\prod_{\left(  i,j\right)  \in
G}\\\text{(since }G=\left\{  \left(  i,j\right)  \in\mathbb{Z}^{2}%
\ \mid\ 1\leq i<j\leq n\right\}  \text{)}}}\underbrace{\left(  x_{j}%
-x_{i}\right)  }_{=\left(  -1\right)  \left(  x_{i}-x_{j}\right)  }\nonumber\\
&  \ \ \ \ \ \ \ \ \ \ \left(  \text{here, we renamed the index }\left(
j,i\right)  \text{ as }\left(  i,j\right)  \text{ in the product}\right)
\nonumber\\
&  =\prod_{\left(  i,j\right)  \in G}\left(  \left(  -1\right)  \left(
x_{i}-x_{j}\right)  \right) \nonumber\\
&  =\underbrace{\left(  -1\right)  ^{\left\vert G\right\vert }}%
_{\substack{=\left(  -1\right)  ^{n\left(  n-1\right)  /2}\\\text{(since
}\left\vert G\right\vert =n\left(  n-1\right)  /2\text{)}}}\underbrace{\prod
_{\left(  i,j\right)  \in G}}_{\substack{=\prod_{\substack{\left(  i,j\right)
\in\mathbb{Z}^{2};\\1\leq i<j\leq n}}\\\text{(since }G=\left\{  \left(
i,j\right)  \in\mathbb{Z}^{2}\ \mid\ 1\leq i<j\leq n\right\}  \text{)}%
}}\left(  x_{i}-x_{j}\right) \nonumber\\
&  =\left(  -1\right)  ^{n\left(  n-1\right)  /2}\underbrace{\prod
_{\substack{\left(  i,j\right)  \in\mathbb{Z}^{2};\\1\leq i<j\leq n}}}%
_{=\prod_{1\leq i<j\leq n}}\left(  x_{i}-x_{j}\right) \nonumber\\
&  =\left(  -1\right)  ^{n\left(  n-1\right)  /2}\prod_{1\leq i<j\leq
n}\left(  x_{i}-x_{j}\right)  . \label{pf.thm.vander-det.d.2}%
\end{align}
Compared with (\ref{pf.thm.vander-det.d.0}), this yields $\det\left(  \left(
x_{j}^{i-1}\right)  _{1\leq i\leq n,\ 1\leq j\leq n}\right)  =\prod_{1\leq
j<i\leq n}\left(  x_{i}-x_{j}\right)  $. This proves Theorem
\ref{thm.vander-det} \textbf{(d)}.
\end{verlong}

\begin{vershort}
\textbf{(c)} We can derive Theorem \ref{thm.vander-det} \textbf{(c)} from
Theorem \ref{thm.vander-det} \textbf{(d)} in the same way as we derived part
\textbf{(b)} from \textbf{(a)}. \qedhere

\end{vershort}

\begin{verlong}
\textbf{(c)} The definition of the transpose of a matrix yields $\left(
\left(  x_{i}^{j-1}\right)  _{1\leq i\leq n,\ 1\leq j\leq n}\right)
^{T}=\left(  x_{j}^{i-1}\right)  _{1\leq i\leq n,\ 1\leq j\leq n}$. Hence,%
\[
\det\left(  \underbrace{\left(  \left(  x_{i}^{j-1}\right)  _{1\leq i\leq
n,\ 1\leq j\leq n}\right)  ^{T}}_{=\left(  x_{j}^{i-1}\right)  _{1\leq i\leq
n,\ 1\leq j\leq n}}\right)  =\det\left(  \left(  x_{j}^{i-1}\right)  _{1\leq
i\leq n,\ 1\leq j\leq n}\right)  =\prod_{1\leq j<i\leq n}\left(  x_{i}%
-x_{j}\right)
\]
(by Theorem \ref{thm.vander-det} \textbf{(d)}). Compared with%
\[
\det\left(  \left(  \left(  x_{i}^{j-1}\right)  _{1\leq i\leq n,\ 1\leq j\leq
n}\right)  ^{T}\right)  =\det\left(  \left(  x_{i}^{j-1}\right)  _{1\leq i\leq
n,\ 1\leq j\leq n}\right)
\]
(by Exercise \ref{exe.ps4.4}, applied to $A=\left(  x_{i}^{j-1}\right)
_{1\leq i\leq n,\ 1\leq j\leq n}$), this yields%
\[
\det\left(  \left(  x_{i}^{j-1}\right)  _{1\leq i\leq n,\ 1\leq j\leq
n}\right)  =\prod_{1\leq j<i\leq n}\left(  x_{i}-x_{j}\right)  .
\]
This proves Theorem \ref{thm.vander-det} \textbf{(c)}.
\end{verlong}
\end{proof}

\subsubsection{\label{subsect.vandermonde.factoring}A proof by factoring the
matrix}

Next, I shall outline another proof of Theorem \ref{thm.vander-det}, which
proceeds by writing the matrix $\left(  x_{j}^{i-1}\right)  _{1\leq i\leq
n,\ 1\leq j\leq n}$ as a product of a lower-triangular matrix with an
upper-triangular matrix. The idea of this proof appears in \cite[Theorem
2.1]{OruPhi}, \cite[Theorem 2]{GohKol} and \cite[proof of Lemma 5.16]{OlvSha}
(although the first two of these three sources use slightly different
arguments, and the third gives only a hint of the proof).

We will need several lemmas for the proof. The proofs of these lemmas are
relegated to the solution of Exercise \ref{exe.subsect.vandermonde.factoring}.

We begin with a definition that will be used throughout Subsection
\ref{subsect.vandermonde.factoring}:

\begin{definition}
\label{def.vandermonde.factoring.h}Let $k\in\mathbb{Z}$ and $n\in\mathbb{N}$.
Let $x_{1},x_{2},\ldots,x_{n}$ be $n$ elements of $\mathbb{K}$. Then, we
define an element $h_{k}\left(  x_{1},x_{2},\ldots,x_{n}\right)  \in
\mathbb{K}$ by%
\[
h_{k}\left(  x_{1},x_{2},\ldots,x_{n}\right)  =\sum_{\substack{\left(
a_{1},a_{2},\ldots,a_{n}\right)  \in\mathbb{N}^{n};\\a_{1}+a_{2}+\cdots
+a_{n}=k}}x_{1}^{a_{1}}x_{2}^{a_{2}}\cdots x_{n}^{a_{n}}.
\]
(Note that the sum on the right hand side of this equality is finite, because
only finitely many $\left(  a_{1},a_{2},\ldots,a_{n}\right)  \in\mathbb{N}%
^{n}$ satisfy $a_{1}+a_{2}+\cdots+a_{n}=k$.)
\end{definition}

The element $h_{k}\left(  x_{1},x_{2},\ldots,x_{n}\right)  $ defined in
Definition \ref{def.vandermonde.factoring.h} is often called the
$k$\textit{-th complete homogeneous function of the }$n$ \textit{elements
}$x_{1},x_{2},\ldots,x_{n}$ (although, more often, this notion is reserved for
a different, more abstract object, of whom $h_{k}\left(  x_{1},x_{2}%
,\ldots,x_{n}\right)  $ is just an evaluation). Let us see some examples:

\begin{example}
\label{exam.vandermonde.factoring.h.12}\textbf{(a)} If $n\in\mathbb{N}$, and
if $x_{1},x_{2},\ldots,x_{n}$ are $n$ elements of $\mathbb{K}$, then%
\[
h_{1}\left(  x_{1},x_{2},\ldots,x_{n}\right)  =x_{1}+x_{2}+\cdots+x_{n}%
\]
and%
\[
h_{2}\left(  x_{1},x_{2},\ldots,x_{n}\right)  =\left(  x_{1}^{2}+x_{2}%
^{2}+\cdots+x_{n}^{2}\right)  +\sum_{1\leq i<j\leq n}x_{i}x_{j}.
\]
For example, for $n=3$, we obtain $h_{2}\left(  x_{1},x_{2},x_{3}\right)
=\left(  x_{1}^{2}+x_{2}^{2}+x_{3}^{2}\right)  +\left(  x_{1}x_{2}+x_{1}%
x_{3}+x_{2}x_{3}\right)  $.

\textbf{(b)} If $x$ and $y$ are two elements of $\mathbb{K}$, then%
\[
h_{k}\left(  x,y\right)  =\sum_{\substack{\left(  a_{1},a_{2}\right)
\in\mathbb{N}^{2};\\a_{1}+a_{2}=k}}x^{a_{1}}y^{a_{2}}=x^{k}y^{0}+x^{k-1}%
y^{1}+\cdots+x^{0}y^{k}%
\]
for every $k\in\mathbb{N}$.

\textbf{(c)} If $x\in\mathbb{K}$, then $h_{k}\left(  x\right)  =x^{k}$ for
every $k\in\mathbb{N}$.
\end{example}

\begin{lemma}
\label{lem.vandermonde.factoring.h.0}Let $n\in\mathbb{N}$. Let $x_{1}%
,x_{2},\ldots,x_{n}$ be $n$ elements of $\mathbb{K}$.

\textbf{(a)} We have $h_{k}\left(  x_{1},x_{2},\ldots,x_{n}\right)  =0$ for
every negative integer $k$.

\textbf{(b)} We have $h_{0}\left(  x_{1},x_{2},\ldots,x_{n}\right)  =1$.
\end{lemma}

As we have said, all lemmas in Subsection \ref{subsect.vandermonde.factoring}
will be proven in the solution to Exercise
\ref{exe.subsect.vandermonde.factoring}.

Three further lemmas on the $h_{k}\left(  x_{1},x_{2},\ldots,x_{n}\right)  $
will be of use:

\begin{lemma}
\label{lem.vandermonde.factoring.h.hq1}Let $k$ be a positive integer. Let
$x_{1},x_{2},\ldots,x_{k}$ be $k$ elements of $\mathbb{K}$. Let $q\in
\mathbb{Z}$. Then,%
\[
h_{q}\left(  x_{1},x_{2},\ldots,x_{k}\right)  =\sum_{r=0}^{q}x_{k}^{r}%
h_{q-r}\left(  x_{1},x_{2},\ldots,x_{k-1}\right)  .
\]

\end{lemma}

\begin{lemma}
\label{lem.vandermonde.factoring.h.hq2}Let $k$ be a positive integer. Let
$x_{1},x_{2},\ldots,x_{k}$ be $k$ elements of $\mathbb{K}$. Let $q\in
\mathbb{Z}$. Then,%
\[
h_{q}\left(  x_{1},x_{2},\ldots,x_{k}\right)  =h_{q}\left(  x_{1},x_{2}%
,\ldots,x_{k-1}\right)  +x_{k}h_{q-1}\left(  x_{1},x_{2},\ldots,x_{k}\right)
.
\]

\end{lemma}

\begin{lemma}
\label{lem.vandermonde.factoring.h.sum-u}Let $i$ be a positive integer. Let
$x_{1},x_{2},\ldots,x_{i}$ be $i$ elements of $\mathbb{K}$. Let $u\in
\mathbb{K}$. Then,%
\[
\sum_{k=1}^{i}h_{i-k}\left(  x_{1},x_{2},\ldots,x_{k}\right)  \prod
_{p=1}^{k-1}\left(  u-x_{p}\right)  =u^{i-1}.
\]

\end{lemma}

Next, let us introduce two matrices:

\begin{lemma}
\label{lem.vandermonde.factoring.U}Let $n\in\mathbb{N}$. Let $x_{1}%
,x_{2},\ldots,x_{n}$ be $n$ elements of $\mathbb{K}$. Define an $n\times
n$-matrix $U\in\mathbb{K}^{n\times n}$ by%
\[
U=\left(  \prod_{p=1}^{i-1}\left(  x_{j}-x_{p}\right)  \right)  _{1\leq i\leq
n,\ 1\leq j\leq n}.
\]
Then, $\det U=\prod_{1\leq j<i\leq n}\left(  x_{i}-x_{j}\right)  $.
\end{lemma}

\begin{example}
If $n=4$, then the matrix $U$ defined in Lemma
\ref{lem.vandermonde.factoring.U} looks as follows:%
\begin{align*}
U  &  =\left(  \prod_{p=1}^{i-1}\left(  x_{j}-x_{p}\right)  \right)  _{1\leq
i\leq4,\ 1\leq j\leq4}\\
&  =\left(
\begin{array}
[c]{cccc}%
\prod_{p=1}^{0}\left(  x_{1}-x_{p}\right)  & \prod_{p=1}^{0}\left(
x_{2}-x_{p}\right)  & \prod_{p=1}^{0}\left(  x_{3}-x_{p}\right)  & \prod
_{p=1}^{0}\left(  x_{4}-x_{p}\right) \\
\prod_{p=1}^{1}\left(  x_{1}-x_{p}\right)  & \prod_{p=1}^{1}\left(
x_{2}-x_{p}\right)  & \prod_{p=1}^{1}\left(  x_{3}-x_{p}\right)  & \prod
_{p=1}^{1}\left(  x_{4}-x_{p}\right) \\
\prod_{p=1}^{2}\left(  x_{1}-x_{p}\right)  & \prod_{p=1}^{2}\left(
x_{2}-x_{p}\right)  & \prod_{p=1}^{2}\left(  x_{3}-x_{p}\right)  & \prod
_{p=1}^{2}\left(  x_{4}-x_{p}\right) \\
\prod_{p=1}^{3}\left(  x_{1}-x_{p}\right)  & \prod_{p=1}^{3}\left(
x_{2}-x_{p}\right)  & \prod_{p=1}^{3}\left(  x_{3}-x_{p}\right)  & \prod
_{p=1}^{3}\left(  x_{4}-x_{p}\right)
\end{array}
\right) \\
&  =\left(
\begin{array}
[c]{cccc}%
1 & 1 & 1 & 1\\
0 & x_{2}-x_{1} & x_{3}-x_{1} & x_{4}-x_{1}\\
0 & 0 & \left(  x_{3}-x_{1}\right)  \left(  x_{3}-x_{2}\right)  & \left(
x_{4}-x_{1}\right)  \left(  x_{4}-x_{2}\right) \\
0 & 0 & 0 & \left(  x_{4}-x_{1}\right)  \left(  x_{4}-x_{2}\right)  \left(
x_{4}-x_{3}\right)
\end{array}
\right)  .
\end{align*}
(Here, we have used the fact that if $i>j$, then the product $\prod
_{p=1}^{i-1}\left(  x_{j}-x_{p}\right)  $ contains the factor $x_{j}-x_{j}=0$
and thus equals $0$.)
\end{example}

\begin{lemma}
\label{lem.vandermonde.factoring.L}Let $n\in\mathbb{N}$. Let $x_{1}%
,x_{2},\ldots,x_{n}$ be $n$ elements of $\mathbb{K}$. Define an $n\times
n$-matrix $L\in\mathbb{K}^{n\times n}$ by%
\[
L=\left(  h_{i-j}\left(  x_{1},x_{2},\ldots,x_{j}\right)  \right)  _{1\leq
i\leq n,\ 1\leq j\leq n}.
\]
Then, $\det L=1$.
\end{lemma}

\begin{example}
If $n=4$, then the matrix $L$ defined in Lemma
\ref{lem.vandermonde.factoring.L} looks as follows:%
\begin{align*}
L  &  =\left(  h_{i-j}\left(  x_{1},x_{2},\ldots,x_{j}\right)  \right)
_{1\leq i\leq n,\ 1\leq j\leq n}\\
&  =\left(
\begin{array}
[c]{cccc}%
h_{0}\left(  x_{1}\right)  & h_{-1}\left(  x_{1},x_{2}\right)  & h_{-2}\left(
x_{1},x_{2},x_{3}\right)  & h_{-3}\left(  x_{1},x_{2},x_{3},x_{4}\right) \\
h_{1}\left(  x_{1}\right)  & h_{0}\left(  x_{1},x_{2}\right)  & h_{-1}\left(
x_{1},x_{2},x_{3}\right)  & h_{-2}\left(  x_{1},x_{2},x_{3},x_{4}\right) \\
h_{2}\left(  x_{2}\right)  & h_{1}\left(  x_{1},x_{2}\right)  & h_{0}\left(
x_{1},x_{2},x_{3}\right)  & h_{-1}\left(  x_{1},x_{2},x_{3},x_{4}\right) \\
h_{3}\left(  x_{3}\right)  & h_{2}\left(  x_{1},x_{2}\right)  & h_{1}\left(
x_{1},x_{2},x_{3}\right)  & h_{0}\left(  x_{1},x_{2},x_{3},x_{4}\right)
\end{array}
\right) \\
&  =\left(
\begin{array}
[c]{cccc}%
1 & 0 & 0 & 0\\
x_{1} & 1 & 0 & 0\\
x_{1}^{2} & x_{1}+x_{2} & 1 & 0\\
x_{1}^{3} & x_{1}^{2}+x_{1}x_{2}+x_{2}^{2} & x_{1}+x_{2}+x_{3} & 1
\end{array}
\right)  .
\end{align*}
(The fact that the diagonal entries are $1$ is a consequence of Lemma
\ref{lem.vandermonde.factoring.h.0} \textbf{(b)}, and the fact that the
entries above the diagonal are $0$ is a consequence of Lemma
\ref{lem.vandermonde.factoring.h.0} \textbf{(a)}.)
\end{example}

\begin{lemma}
\label{lem.vandermonde.factoring.A=LU}Let $n\in\mathbb{N}$. Let $x_{1}%
,x_{2},\ldots,x_{n}$ be $n$ elements of $\mathbb{K}$. Let $L$ be the $n\times
n$-matrix defined in Lemma \ref{lem.vandermonde.factoring.L}. Let $U$ be the
$n\times n$-matrix defined in Lemma \ref{lem.vandermonde.factoring.U}. Then,%
\[
\left(  x_{j}^{i-1}\right)  _{1\leq i\leq n,\ 1\leq j\leq n}=LU.
\]

\end{lemma}

\begin{exercise}
\label{exe.subsect.vandermonde.factoring}Prove Lemma
\ref{lem.vandermonde.factoring.h.0}, Lemma
\ref{lem.vandermonde.factoring.h.hq1}, Lemma
\ref{lem.vandermonde.factoring.h.hq2}, Lemma
\ref{lem.vandermonde.factoring.h.sum-u}, Lemma
\ref{lem.vandermonde.factoring.U}, Lemma \ref{lem.vandermonde.factoring.L} and
Lemma \ref{lem.vandermonde.factoring.A=LU}.
\end{exercise}

Now, we can prove Theorem \ref{thm.vander-det} again:

\begin{proof}
[Second proof of Theorem \ref{thm.vander-det}.]\textbf{(d)} Let $L$ be the
$n\times n$-matrix defined in Lemma \ref{lem.vandermonde.factoring.L}. Let $U$
be the $n\times n$-matrix defined in Lemma \ref{lem.vandermonde.factoring.U}.
Then, Lemma \ref{lem.vandermonde.factoring.A=LU} yields $\left(  x_{j}%
^{i-1}\right)  _{1\leq i\leq n,\ 1\leq j\leq n}=LU$. Hence,%
\begin{align*}
\det\left(  \underbrace{\left(  x_{j}^{i-1}\right)  _{1\leq i\leq n,\ 1\leq
j\leq n}}_{=LU}\right)   &  =\det\left(  LU\right)  =\underbrace{\det
L}_{\substack{=1\\\text{(by Lemma \ref{lem.vandermonde.factoring.L})}}%
}\cdot\underbrace{\det U}_{\substack{=\prod_{1\leq j<i\leq n}\left(
x_{i}-x_{j}\right)  \\\text{(by Lemma \ref{lem.vandermonde.factoring.U})}}}\\
&  \ \ \ \ \ \ \ \ \ \ \left(
\begin{array}
[c]{c}%
\text{by Theorem \ref{thm.det(AB)}, applied to }L\text{ and }U\\
\text{instead of }A\text{ and }B
\end{array}
\right) \\
&  =1\cdot\prod_{1\leq j<i\leq n}\left(  x_{i}-x_{j}\right)  =\prod_{1\leq
j<i\leq n}\left(  x_{i}-x_{j}\right)  .
\end{align*}
This proves Theorem \ref{thm.vander-det} \textbf{(d)}.

\begin{vershort}
Now, it remains to prove parts \textbf{(a)}, \textbf{(b)} and \textbf{(c)} of
Theorem \ref{thm.vander-det}. This is fairly easy: Back in our First proof of
Theorem \ref{thm.vander-det}, we have derived parts \textbf{(b)}, \textbf{(d)}
and \textbf{(c)} from part \textbf{(a)}. By essentially the same arguments
(sometimes done in reverse), we can derive parts \textbf{(a)}, \textbf{(b)}
and \textbf{(c)} from part \textbf{(d)}. (We need to use the fact that
$\left(  \left(  -1\right)  ^{n\left(  n-1\right)  /2}\right)  ^{2}=1$.)
\qedhere

\end{vershort}

\begin{verlong}
\textbf{(b)} Applying Lemma \ref{lem.vander-det.lem-rearr} to $a_{i,j}%
=x_{j}^{i-1}$, we obtain%
\begin{align*}
\det\left(  \left(  x_{j}^{\left(  n+1-i\right)  -1}\right)  _{1\leq i\leq
n,\ 1\leq j\leq n}\right)   &  =\left(  -1\right)  ^{n\left(  n-1\right)
/2}\underbrace{\det\left(  \left(  x_{j}^{i-1}\right)  _{1\leq i\leq n,\ 1\leq
j\leq n}\right)  }_{\substack{=\prod_{1\leq j<i\leq n}\left(  x_{i}%
-x_{j}\right)  \\\text{(by Theorem \ref{thm.vander-det} \textbf{(d)})}}}\\
&  =\left(  -1\right)  ^{n\left(  n-1\right)  /2}\prod_{1\leq j<i\leq
n}\left(  x_{i}-x_{j}\right)  .
\end{align*}
This rewrites as%
\begin{equation}
\det\left(  \left(  x_{j}^{n-i}\right)  _{1\leq i\leq n,\ 1\leq j\leq
n}\right)  =\left(  -1\right)  ^{n\left(  n-1\right)  /2}\prod_{1\leq j<i\leq
n}\left(  x_{i}-x_{j}\right)  \label{pf.thm.vander-det.pf2.b.0}%
\end{equation}
(since every $\left(  i,j\right)  \in\left\{  1,2,\ldots,n\right\}  ^{2}$
satisfies $x_{j}^{\left(  n+1-i\right)  -1}=x_{j}^{n-i}$).

But we have%
\begin{equation}
\prod_{1\leq j<i\leq n}\left(  x_{i}-x_{j}\right)  =\left(  -1\right)
^{n\left(  n-1\right)  /2}\prod_{1\leq i<j\leq n}\left(  x_{i}-x_{j}\right)  .
\label{pf.thm.vander-det.pf2.b.1}%
\end{equation}
(Indeed, this is the equality (\ref{pf.thm.vander-det.d.2}) from the First
proof of Theorem \ref{thm.vander-det} above.) Also, the number $n\left(
n-1\right)  /2$ is an integer. Hence, $n\left(  n-1\right)  $ is even. Thus,
$\left(  -1\right)  ^{n\left(  n-1\right)  }=1$. Now,%
\begin{align*}
\left(  \left(  -1\right)  ^{n\left(  n-1\right)  /2}\right)  ^{2}  &
=\left(  -1\right)  ^{\left(  n\left(  n-1\right)  /2\right)  \cdot2}=\left(
-1\right)  ^{n\left(  n-1\right)  }\ \ \ \ \ \ \ \ \ \ \left(  \text{since
}\left(  n\left(  n-1\right)  /2\right)  \cdot2=n\left(  n-1\right)  \right)
\\
&  =1.
\end{align*}
Now, (\ref{pf.thm.vander-det.pf2.b.0}) becomes%
\begin{align*}
\det\left(  \left(  x_{j}^{n-i}\right)  _{1\leq i\leq n,\ 1\leq j\leq
n}\right)   &  =\left(  -1\right)  ^{n\left(  n-1\right)  /2}\underbrace{\prod
_{1\leq j<i\leq n}\left(  x_{i}-x_{j}\right)  }_{\substack{=\left(  -1\right)
^{n\left(  n-1\right)  /2}\prod_{1\leq i<j\leq n}\left(  x_{i}-x_{j}\right)
\\\text{(by (\ref{pf.thm.vander-det.pf2.b.1}))}}}\\
&  =\underbrace{\left(  -1\right)  ^{n\left(  n-1\right)  /2}\left(
-1\right)  ^{n\left(  n-1\right)  /2}}_{=\left(  \left(  -1\right)  ^{n\left(
n-1\right)  /2}\right)  ^{2}=1}\prod_{1\leq i<j\leq n}\left(  x_{i}%
-x_{j}\right) \\
&  =\prod_{1\leq i<j\leq n}\left(  x_{i}-x_{j}\right)  .
\end{align*}
This proves Theorem \ref{thm.vander-det} \textbf{(b)}.

\textbf{(a)} The definition of the transpose of a matrix yields $\left(
\left(  x_{j}^{n-i}\right)  _{1\leq i\leq n,\ 1\leq j\leq n}\right)
^{T}=\left(  x_{i}^{n-j}\right)  _{1\leq i\leq n,\ 1\leq j\leq n}$. Hence,%
\[
\det\left(  \underbrace{\left(  \left(  x_{j}^{n-i}\right)  _{1\leq i\leq
n,\ 1\leq j\leq n}\right)  ^{T}}_{=\left(  x_{i}^{n-j}\right)  _{1\leq i\leq
n,\ 1\leq j\leq n}}\right)  =\det\left(  \left(  x_{i}^{n-j}\right)  _{1\leq
i\leq n,\ 1\leq j\leq n}\right)  .
\]
Hence,%
\begin{align*}
\det\left(  \left(  x_{i}^{n-j}\right)  _{1\leq i\leq n,\ 1\leq j\leq
n}\right)   &  =\det\left(  \left(  \left(  x_{j}^{n-i}\right)  _{1\leq i\leq
n,\ 1\leq j\leq n}\right)  ^{T}\right)  =\det\left(  \left(  x_{j}%
^{n-i}\right)  _{1\leq i\leq n,\ 1\leq j\leq n}\right) \\
&  \ \ \ \ \ \ \ \ \ \ \left(  \text{by Exercise \ref{exe.ps4.4}, applied to
}A=\left(  x_{j}^{n-i}\right)  _{1\leq i\leq n,\ 1\leq j\leq n}\right) \\
&  =\prod_{1\leq i<j\leq n}\left(  x_{i}-x_{j}\right)  .
\end{align*}
This proves Theorem \ref{thm.vander-det} \textbf{(a)}.

\textbf{(c)} The definition of the transpose of a matrix yields $\left(
\left(  x_{i}^{j-1}\right)  _{1\leq i\leq n,\ 1\leq j\leq n}\right)
^{T}=\left(  x_{j}^{i-1}\right)  _{1\leq i\leq n,\ 1\leq j\leq n}$. Hence,%
\[
\det\left(  \underbrace{\left(  \left(  x_{i}^{j-1}\right)  _{1\leq i\leq
n,\ 1\leq j\leq n}\right)  ^{T}}_{=\left(  x_{j}^{i-1}\right)  _{1\leq i\leq
n,\ 1\leq j\leq n}}\right)  =\det\left(  \left(  x_{j}^{i-1}\right)  _{1\leq
i\leq n,\ 1\leq j\leq n}\right)  =\prod_{1\leq j<i\leq n}\left(  x_{i}%
-x_{j}\right)
\]
(by Theorem \ref{thm.vander-det} \textbf{(d)}). Compared with%
\[
\det\left(  \left(  \left(  x_{i}^{j-1}\right)  _{1\leq i\leq n,\ 1\leq j\leq
n}\right)  ^{T}\right)  =\det\left(  \left(  x_{i}^{j-1}\right)  _{1\leq i\leq
n,\ 1\leq j\leq n}\right)
\]
(by Exercise \ref{exe.ps4.4}, applied to $A=\left(  x_{i}^{j-1}\right)
_{1\leq i\leq n,\ 1\leq j\leq n}$), this yields%
\[
\det\left(  \left(  x_{i}^{j-1}\right)  _{1\leq i\leq n,\ 1\leq j\leq
n}\right)  =\prod_{1\leq j<i\leq n}\left(  x_{i}-x_{j}\right)  .
\]
This proves Theorem \ref{thm.vander-det} \textbf{(c)}.
\end{verlong}
\end{proof}

\subsubsection{Remarks and variations}

\begin{remark}
\label{rmk.vander-det.sign}One consequence of Theorem \ref{thm.vander-det} is
a new solution to Exercise \ref{exe.perm.sign.pseudoexplicit} \textbf{(a)}:

Namely, let $n\in\mathbb{N}$ and $\sigma\in S_{n}$. Let $x_{1},x_{2}%
,\ldots,x_{n}$ be $n$ elements of $\mathbb{C}$ (or of any commutative ring).
Then, Theorem \ref{thm.vander-det} \textbf{(a)} yields
\[
\det\left(  \left(  x_{i}^{n-j}\right)  _{1\leq i\leq n,\ 1\leq j\leq
n}\right)  =\prod_{1\leq i<j\leq n}\left(  x_{i}-x_{j}\right)  .
\]
On the other hand, Theorem \ref{thm.vander-det} \textbf{(a)} (applied to
$x_{\sigma\left(  1\right)  },x_{\sigma\left(  2\right)  },\ldots
,x_{\sigma\left(  n\right)  }$ instead of $x_{1},x_{2},\ldots,x_{n}$) yields%
\begin{equation}
\det\left(  \left(  x_{\sigma\left(  i\right)  }^{n-j}\right)  _{1\leq i\leq
n,\ 1\leq j\leq n}\right)  =\prod_{1\leq i<j\leq n}\left(  x_{\sigma\left(
i\right)  }-x_{\sigma\left(  j\right)  }\right)  .
\label{eq.rmk.vander-det.sign.2}%
\end{equation}
But Lemma \ref{lem.det.sigma} \textbf{(a)} (applied to $B=\left(  x_{i}%
^{n-j}\right)  _{1\leq i\leq n,\ 1\leq j\leq n}$, $\kappa=\sigma$ and
$B_{\kappa}=\left(  x_{\sigma\left(  i\right)  }^{n-j}\right)  _{1\leq i\leq
n,\ 1\leq j\leq n}$) yields%
\begin{align*}
\det\left(  \left(  x_{\sigma\left(  i\right)  }^{n-j}\right)  _{1\leq i\leq
n,\ 1\leq j\leq n}\right)   &  =\left(  -1\right)  ^{\sigma}\cdot
\underbrace{\det\left(  \left(  x_{i}^{n-j}\right)  _{1\leq i\leq n,\ 1\leq
j\leq n}\right)  }_{=\prod_{1\leq i<j\leq n}\left(  x_{i}-x_{j}\right)  }\\
&  =\left(  -1\right)  ^{\sigma}\cdot\prod_{1\leq i<j\leq n}\left(
x_{i}-x_{j}\right)  .
\end{align*}
Compared with (\ref{eq.rmk.vander-det.sign.2}), this yields%
\[
\prod_{1\leq i<j\leq n}\left(  x_{\sigma\left(  i\right)  }-x_{\sigma\left(
j\right)  }\right)  =\left(  -1\right)  ^{\sigma}\cdot\prod_{1\leq i<j\leq
n}\left(  x_{i}-x_{j}\right)  .
\]
Thus, Exercise \ref{exe.perm.sign.pseudoexplicit} \textbf{(a)} is solved.
However, Exercise \ref{exe.perm.sign.pseudoexplicit} \textbf{(b)} cannot be
solved this way.
\end{remark}

\begin{exercise}
\label{exe.vander-det.s1}Let $n$ be a positive integer. Let $x_{1}%
,x_{2},\ldots,x_{n}$ be $n$ elements of $\mathbb{K}$. Prove that%
\[
\det\left(  \left(
\begin{cases}
x_{i}^{n-j}, & \text{if }j>1;\\
x_{i}^{n}, & \text{if }j=1
\end{cases}
\right)  _{1\leq i\leq n,\ 1\leq j\leq n}\right)  =\left(  x_{1}+x_{2}%
+\cdots+x_{n}\right)  \prod_{1\leq i<j\leq n}\left(  x_{i}-x_{j}\right)  .
\]
(For example, when $n=4$, this states that%
\begin{align*}
&  \det\left(
\begin{array}
[c]{cccc}%
x_{1}^{4} & x_{1}^{2} & x_{1} & 1\\
x_{2}^{4} & x_{2}^{2} & x_{2} & 1\\
x_{3}^{4} & x_{3}^{2} & x_{3} & 1\\
x_{4}^{4} & x_{4}^{2} & x_{4} & 1
\end{array}
\right) \\
&  =\left(  x_{1}+x_{2}+x_{3}+x_{4}\right)  \left(  x_{1}-x_{2}\right)
\left(  x_{1}-x_{3}\right)  \left(  x_{1}-x_{4}\right)  \left(  x_{2}%
-x_{3}\right)  \left(  x_{2}-x_{4}\right)  \left(  x_{3}-x_{4}\right)  .
\end{align*}
)
\end{exercise}

\begin{remark}
\label{rmk.vander-det.schur}We can try to generalize Vandermonde's
determinant. Namely, let $n\in\mathbb{N}$. Let $x_{1},x_{2},\ldots,x_{n}$ be
$n$ elements of $\mathbb{K}$. Let $a_{1},a_{2},\ldots,a_{n}$ be $n$
nonnegative integers. Let $A$ be the $n\times n$-matrix
\[
\left(  x_{i}^{a_{j}}\right)  _{1\leq i\leq n,\ 1\leq j\leq n}=\left(
\begin{array}
[c]{cccc}%
x_{1}^{a_{1}} & x_{1}^{a_{2}} & \cdots & x_{1}^{a_{n}}\\
x_{2}^{a_{1}} & x_{2}^{a_{2}} & \cdots & x_{2}^{a_{n}}\\
\vdots & \vdots & \ddots & \vdots\\
x_{n}^{a_{1}} & x_{n}^{a_{2}} & \cdots & x_{n}^{a_{n}}%
\end{array}
\right)  .
\]
What can we say about $\det A$ ?

Theorem \ref{thm.vander-det} says that if $\left(  a_{1},a_{2},\ldots
,a_{n}\right)  =\left(  n-1,n-2,\ldots,0\right)  $, then $\det A=\prod_{1\leq
i<j\leq n}\left(  x_{i}-x_{j}\right)  $.

Exercise \ref{exe.vander-det.s1} says that if $n>0$ and $\left(  a_{1}%
,a_{2},\ldots,a_{n}\right)  =\left(  n,n-2,n-3,\ldots,0\right)  $, then $\det
A=\left(  x_{1}+x_{2}+\cdots+x_{n}\right)  \prod_{1\leq i<j\leq n}\left(
x_{i}-x_{j}\right)  $.

This suggests a general pattern: We would suspect that for every $\left(
a_{1},a_{2},\ldots,a_{n}\right)  $, there is a polynomial $P_{\left(
a_{1},a_{2},\ldots,a_{n}\right)  }$ in $n$ indeterminates $X_{1},X_{2}%
,\ldots,X_{n}$ such that%
\[
\det A=P_{\left(  a_{1},a_{2},\ldots,a_{n}\right)  }\left(  x_{1},x_{2}%
,\ldots,x_{n}\right)  \cdot\prod_{1\leq i<j\leq n}\left(  x_{i}-x_{j}\right)
.
\]

It turns out that this is true. Moreover, this polynomial $P_{\left(
a_{1},a_{2},\ldots,a_{n}\right)  }$ is:

\begin{itemize}
\item zero if two of $a_{1},a_{2},\ldots,a_{n}$ are equal;

\item homogeneous of degree $a_{1}+a_{2}+\cdots+a_{n}-\dbinom{n}{2}$;

\item symmetric in $X_{1},X_{2},\ldots,X_{n}$.
\end{itemize}

For example,%
\begin{align*}
P_{\left(  n-1,n-2,\ldots,0\right)  }  &  =1;\\
P_{\left(  n,n-2,n-3,\ldots,0\right)  }  &  =\sum_{i=1}^{n}X_{i}=X_{1}%
+X_{2}+\cdots+X_{n};\\
P_{\left(  n,n-1,\ldots,n-k+1,n-k-1,n-k-2,\ldots,0\right)  }  &  =\sum_{1\leq
i_{1}<i_{2}<\cdots<i_{k}\leq n}X_{i_{1}}X_{i_{2}}\cdots X_{i_{k}}\\
&  \ \ \ \ \ \ \ \ \ \ \text{for every }k\in\left\{  0,1,\ldots,n\right\}  ;\\
P_{\left(  n+1,n-2,n-3,\ldots,0\right)  }  &  =\sum_{1\leq i\leq j\leq n}%
X_{i}X_{j};\\
P_{\left(  n+1,n-1,n-3,n-4,\ldots,0\right)  }  &  =\sum_{1\leq i<j\leq
n}\left(  X_{i}^{2}X_{j}+X_{i}X_{j}^{2}\right)  +2\sum_{1\leq i<j<k\leq
n}X_{i}X_{j}X_{k}.
\end{align*}

But this polynomial $P_{\left(  a_{1},a_{2},\ldots,a_{n}\right)  }$ can
actually be described rather explicitly for general $\left(  a_{1}%
,a_{2},\ldots,a_{n}\right)  $; it is a so-called \textit{Schur polynomial} (at
least when $a_{1}>a_{2}>\cdots>a_{n}$; otherwise it is either zero or $\pm$ a
Schur polynomial). See \cite[The Bi-Alternant Formula]{Stembridge},
\cite[Theorem 7.15.1]{Stanley-EC2}, \cite{Leeuwen-aS} or \cite[Theorem
7.3.11]{21s} for the details. (Notice that \cite{Leeuwen-aS} uses the notation
$\varepsilon\left(  \sigma\right)  $ for the sign of a permutation $\sigma$.)
The theory of Schur polynomials shows, in particular, that all coefficients of
the polynomial $P_{\left(  a_{1},a_{2},\ldots,a_{n}\right)  }$ have equal sign
(which is positive if $a_{1}>a_{2}>\cdots>a_{n}$).
\end{remark}

\begin{remark}
\label{rmk.vander-det.secret-ints}There are plenty other variations on the
Vandermonde determinant. For instance, one can try replacing the powers
$x_{i}^{j-1}$ by binomial coefficients $\dbinom{x_{i}}{j-1}$ in Theorem
\ref{thm.vander-det} \textbf{(c)}, at least when these binomial coefficients
are well-defined (e.g., when the $x_{1},x_{2},\ldots,x_{n}$ are complex
numbers). The result is rather nice: If $x_{1},x_{2},\ldots,x_{n}$ are any $n$
complex numbers, then%
\[
\prod_{1\leq i<j\leq n}\dfrac{x_{i}-x_{j}}{i-j}=\det\left(  \left(
\dbinom{x_{i}}{j-1}\right)  _{1\leq i\leq n,\ 1\leq j\leq n}\right)  .
\]
(This equality is proved, e.g., in \cite[Corollary 11]{GriHyp} and in
\cite[\S 9, Example 5]{AndDos}, and follows easily from Exercise
\ref{exe.det.vdm-pol}; it also appeared in \cite[exercise 269]{FadSom72}.)
This equality has the surprising consequence that, whenever $x_{1}%
,x_{2},\ldots,x_{n}$ are $n$ integers, the product $\prod_{1\leq i<j\leq
n}\dfrac{x_{i}-x_{j}}{i-j}$ is an integer as well (because it is the
determinant of a matrix whose entries are integers). This is a nontrivial
result! (A more elementary proof appears in \cite[\S 3, Example 8]{AndDos}.
See also \cite{Sury95} for a proof using cyclotomic polynomials, and
\cite[Theorem 3]{Bharga00} for a placement of this result in a more general context.)

Another \textquotedblleft secret integer\textquotedblright\ (i.e., rational
number which turns out to be an integer for non-obvious reasons) is%
\begin{equation}
\dfrac{H\left(  a\right)  H\left(  b\right)  H\left(  c\right)  H\left(
a+b+c\right)  }{H\left(  b+c\right)  H\left(  c+a\right)  H\left(  a+b\right)
}, \label{eq.rmk.vander-det.secret-ints.H}%
\end{equation}
where $a,b,c$ are three nonnegative integers, and where $H\left(  n\right)  $
(for $n\in\mathbb{N}$) denotes the \textit{hyperfactorial} of $n$, defined by%
\[
H\left(  n\right)  =\prod_{k=0}^{n-1}k!=0!\cdot1!\cdot\cdots\cdot\left(
n-1\right)  !.
\]
I am aware of two proofs of the fact that
(\ref{eq.rmk.vander-det.secret-ints.H}) gives an integer for every
$a,b,c\in\mathbb{N}$: One proof is combinatorial, and argues that
(\ref{eq.rmk.vander-det.secret-ints.H}) is the number of \textit{plane
partitions inside an }$a\times b\times c$\textit{-box} (see \cite[last
equality in \S 7.21]{Stanley-EC2} for a proof), or, equivalently, the number
of \textit{rhombus tilings of a hexagon with sidelengths }$a,b,c,a,b,c$ (see
\cite{Eisenk-planepartits} for a precise statement). Another proof (see
\cite[Theorem 0]{GriHyp}) exhibits (\ref{eq.rmk.vander-det.secret-ints.H}) as
the determinant of a matrix, again using the Vandermonde determinant!
\end{remark}

For some more exercises related to Vandermonde determinants, see \cite[Chapter
1, problems 1.12--1.22]{Prasolov}. Here comes one of them (\cite[Lemma 9 and
Lemma 10]{Kratte05}):

\begin{exercise}
\label{exe.vander-det.xi+yj}Let $n$ be a positive integer. Let $x_{1}%
,x_{2},\ldots,x_{n}$ be $n$ elements of $\mathbb{K}$. Let $y_{1},y_{2}%
,\ldots,y_{n}$ be $n$ elements of $\mathbb{K}$.

\textbf{(a)} For every $m\in\left\{  0,1,\ldots,n-2\right\}  $, prove that%
\[
\det\left(  \left(  \left(  x_{i}+y_{j}\right)  ^{m}\right)  _{1\leq i\leq
n,\ 1\leq j\leq n}\right)  =0.
\]

\textbf{(b)} Prove that%
\begin{align*}
&  \det\left(  \left(  \left(  x_{i}+y_{j}\right)  ^{n-1}\right)  _{1\leq
i\leq n,\ 1\leq j\leq n}\right) \\
&  =\left(  \prod_{k=0}^{n-1}\dbinom{n-1}{k}\right)  \left(  \prod_{1\leq
i<j\leq n}\left(  x_{i}-x_{j}\right)  \right)  \left(  \prod_{1\leq i<j\leq
n}\left(  y_{j}-y_{i}\right)  \right)  .
\end{align*}

[\textbf{Hint:} Use the binomial theorem.]

\textbf{(c)} Let $\left(  p_{0},p_{1},\ldots,p_{n-1}\right)  \in\mathbb{K}%
^{n}$ be an $n$-tuple of elements of $\mathbb{K}$. Let $P\left(  X\right)
\in\mathbb{K}\left[  X\right]  $ be the polynomial $\sum_{k=0}^{n-1}p_{k}%
X^{k}$. Prove that%
\begin{align*}
&  \det\left(  \left(  P\left(  x_{i}+y_{j}\right)  \right)  _{1\leq i\leq
n,\ 1\leq j\leq n}\right) \\
&  =p_{n-1}^{n}\left(  \prod_{k=0}^{n-1}\dbinom{n-1}{k}\right)  \left(
\prod_{1\leq i<j\leq n}\left(  x_{i}-x_{j}\right)  \right)  \left(
\prod_{1\leq i<j\leq n}\left(  y_{j}-y_{i}\right)  \right)  .
\end{align*}

\textbf{(d)} Let $\left(  p_{0},p_{1},\ldots,p_{n-1}\right)  \in\mathbb{K}%
^{n}$ be an $n$-tuple of elements of $\mathbb{K}$. Let $P\left(  X\right)
\in\mathbb{K}\left[  X\right]  $ be the polynomial $\sum_{k=0}^{n-1}p_{k}%
X^{k}$. Prove that%
\begin{align*}
&  \det\left(  \left(  P\left(  x_{i}y_{j}\right)  \right)  _{1\leq i\leq
n,\ 1\leq j\leq n}\right) \\
&  =\left(  \prod_{k=0}^{n-1}p_{k}\right)  \left(  \prod_{1\leq i<j\leq
n}\left(  x_{i}-x_{j}\right)  \right)  \left(  \prod_{1\leq i<j\leq n}\left(
y_{i}-y_{j}\right)  \right)  .
\end{align*}

\end{exercise}

Notice how Exercise \ref{exe.vander-det.xi+yj} \textbf{(a)} generalizes
Example \ref{exam.xi+yj} (for $n\geq3$).

\subsection{\label{sect.invertible}Invertible elements in commutative rings,
and fields}

We shall now interrupt our study of determinants for a moment. Let us define
the notion of inverses in $\mathbb{K}$. (Recall that $\mathbb{K}$ is a
commutative ring.)

\begin{definition}
\label{def.rings.inverse}Let $a\in\mathbb{K}$. Then, an element $b\in
\mathbb{K}$ is said to be an \textit{inverse} of $a$ if it satisfies $ab=1$
and $ba=1$.
\end{definition}

Of course, the two conditions $ab=1$ and $ba=1$ in Definition
\ref{def.rings.inverse} are equivalent, since $ab=ba$ for every $a\in
\mathbb{K}$ and $b\in\mathbb{K}$. Nevertheless, we have given both conditions,
because this way the similarity between the inverse of an element of
$\mathbb{K}$ and the inverse of a map becomes particularly clear.

For example, the element $1$ of $\mathbb{Z}$ is its own inverse (since
$1\cdot1=1$), and the element $-1$ of $\mathbb{Z}$ is its own inverse as well
(since $\left(  -1\right)  \cdot\left(  -1\right)  =1$). These elements $1$
and $-1$ are the only elements of $\mathbb{Z}$ which have an inverse in
$\mathbb{Z}$. However, in the larger commutative ring $\mathbb{Q}$, every
nonzero element $a$ has an inverse (namely, $\dfrac{1}{a}$).

\begin{proposition}
\label{prop.rings.inverse-uni}Let $a\in\mathbb{K}$. Then, there exists at most
one inverse of $a$ in $\mathbb{K}$.
\end{proposition}

\begin{proof}
[Proof of Proposition \ref{prop.rings.inverse-uni}.]Let $b$ and $b^{\prime}$
be any two inverses of $a$ in $\mathbb{K}$. Since $b$ is an inverse of $a$ in
$\mathbb{K}$, we have $ab=1$ and $ba=1$ (by the definition of an
\textquotedblleft inverse of $a$\textquotedblright). Since $b^{\prime}$ is an
inverse of $a$ in $\mathbb{K}$, we have $ab^{\prime}=1$ and $b^{\prime}a=1$
(by the definition of an \textquotedblleft inverse of $a$\textquotedblright).
Now, comparing $b\underbrace{ab^{\prime}}_{=1}=b$ with $\underbrace{ba}%
_{=1}b^{\prime}=b^{\prime}$, we obtain $b=b^{\prime}$.

Let us now forget that we fixed $b$ and $b^{\prime}$. We thus have shown that
if $b$ and $b^{\prime}$ are two inverses of $a$ in $\mathbb{K}$, then
$b=b^{\prime}$. In other words, any two inverses of $a$ in $\mathbb{K}$ are
equal. In other words, there exists at most one inverse of $a$ in $\mathbb{K}%
$. This proves Proposition \ref{prop.rings.inverse-uni}.
\end{proof}

\begin{definition}
\label{def.rings.invertible}\textbf{(a)} An element $a\in\mathbb{K}$ is said
to be \textit{invertible} (or, more precisely, \textit{invertible in
}$\mathbb{K}$) if and only if there exists an inverse of $a$ in $\mathbb{K}$.
In this case, this inverse of $a$ is unique (by Proposition
\ref{prop.rings.inverse-uni}), and thus will be called \textit{the inverse of
}$a$ and denoted by $a^{-1}$.

\textbf{(b)} It is clear that the unity $1$ of $\mathbb{K}$ is invertible
(having inverse $1$). Also, the product of any two invertible elements $a$ and
$b$ of $\mathbb{K}$ is again invertible (having inverse $\left(  ab\right)
^{-1}=a^{-1}b^{-1}$).

\textbf{(c)} If $a$ and $b$ are two elements of $\mathbb{K}$ such that $a$ is
invertible (in $\mathbb{K}$), then we write $\dfrac{b}{a}$ (or $b/a$) for the
product $ba^{-1}$. These fractions behave just as fractions of integers
behave: For example, if $a,b,c,d$ are four elements of $\mathbb{K}$ such that
$a$ and $c$ are invertible, then $\dfrac{b}{a}+\dfrac{d}{c}=\dfrac{bc+da}{ac}$
and $\dfrac{b}{a}\cdot\dfrac{d}{c}=\dfrac{bd}{ac}$ (and the product $ac$ is
indeed invertible, so that the fractions $\dfrac{bc+da}{ac}$ and $\dfrac
{bd}{ac}$ actually make sense).
\end{definition}

Of course, the meaning of the word \textquotedblleft
invertible\textquotedblright\ depends on the ring $\mathbb{K}$. For example,
the integer $2$ is invertible in $\mathbb{Q}$ (because $\dfrac{1}{2}$ is an
inverse of $2$ in $\mathbb{Q}$), but not invertible in $\mathbb{Z}$ (since it
has no inverse in $\mathbb{Z}$). Thus, it is important to say
\textquotedblleft invertible in $\mathbb{K}$\textquotedblright\ unless the
context makes it clear what $\mathbb{K}$ is.

One can usually work with invertible elements in commutative rings in the same
way as one works with nonzero rational numbers. For example, if $a$ is an
invertible element of $\mathbb{K}$, then we can define $a^{n}$ not only for
all $n\in\mathbb{N}$, but also for all $n\in\mathbb{Z}$ (by setting
$a^{n}=\left(  a^{-1}\right)  ^{-n}$ for all negative integers $n$). Of
course, when $n=-1$, this is consistent with our notation $a^{-1}$ for the
inverse of $a$.

Next, we define the notion of a \textit{field}\footnote{We are going to use
the following simple fact: A commutative ring $\mathbb{K}$ is a trivial ring
if and only if $0_{\mathbb{K}}=1_{\mathbb{K}}$.
\par
\textit{Proof.} Assume that $\mathbb{K}$ is a trivial ring. Thus, $\mathbb{K}$
has only one element. Hence, both $0_{\mathbb{K}}$ and $1_{\mathbb{K}}$ have
to equal this one element. Therefore, $0_{\mathbb{K}}=1_{\mathbb{K}}$.
\par
Now, forget that we assumed that $\mathbb{K}$ is a trivial ring. We thus have
proven that
\begin{equation}
\text{if }\mathbb{K}\text{ is a trivial ring, then }0_{\mathbb{K}%
}=1_{\mathbb{K}}\text{.} \label{eq.fn.field.trivring.1}%
\end{equation}
\par
Conversely, assume that $0_{\mathbb{K}}=1_{\mathbb{K}}$. Then, every
$a\in\mathbb{K}$ satisfies $a=a\cdot\underbrace{1_{\mathbb{K}}}%
_{=0_{\mathbb{K}}}=a\cdot0_{\mathbb{K}}=0_{\mathbb{K}}\in\left\{
0_{\mathbb{K}}\right\}  $. In other words, $\mathbb{K}\subseteq\left\{
0_{\mathbb{K}}\right\}  $. Combining this with $\left\{  0_{\mathbb{K}%
}\right\}  \subseteq\mathbb{K}$, we obtain $\mathbb{K}=\left\{  0_{\mathbb{K}%
}\right\}  $. Hence, $\mathbb{K}$ has only one element. In other words,
$\mathbb{K}$ is a trivial ring.
\par
Now, forget that we assumed that $0_{\mathbb{K}}=1_{\mathbb{K}}$. We thus have
proven that
\[
\text{if }0_{\mathbb{K}}=1_{\mathbb{K}}\text{, then }\mathbb{K}\text{ is a
trivial ring.}%
\]
Combining this with (\ref{eq.fn.field.trivring.1}), we conclude that
$\mathbb{K}$ is a trivial ring if and only if $0_{\mathbb{K}}=1_{\mathbb{K}}%
$.}.

\begin{definition}
A commutative ring $\mathbb{K}$ is said to be a \textit{field} if it satisfies
the following two properties:

\begin{itemize}
\item We have $0_{\mathbb{K}}\neq1_{\mathbb{K}}$ (that is, $\mathbb{K}$ is not
a trivial ring).

\item Every element of $\mathbb{K}$ is either zero or invertible.
\end{itemize}
\end{definition}

For example, $\mathbb{Q}$, $\mathbb{R}$ and $\mathbb{C}$ are fields, whereas
polynomial rings such as $\mathbb{Q}\left[  x\right]  $ or $\mathbb{R}\left[
a,b\right]  $ are not fields\footnote{For example, the polynomial $x$ is not
invertible in $\mathbb{Q}\left[  x\right]  $.}. For $n$ being a positive
integer, the ring $\mathbb{Z}/n\mathbb{Z}$ (that is, the ring of residue
classes of integers modulo $n$) is a field if and only if $n$ is a prime number.

Linear algebra (i.e., the study of matrices and linear transformations)
becomes much easier (in many aspects) when $\mathbb{K}$ is a
field\footnote{Many properties of a matrix over a field (such as its rank) are
not even well-defined over an arbitrary commutative ring.}. This is one of the
main reasons why most courses on linear algebra work over fields only (or
begin by working over fields and only later move to the generality of
commutative rings). In these notes we are almost completely limiting ourselves
to the parts of matrix theory which work over any commutative ring.
Nevertheless, let us comment on how determinants can be computed fast when
$\mathbb{K}$ is a field.

\Needspace{30\baselineskip}

\begin{remark}
\label{rmk.laplace.pre.gauss-alg}Assume that $\mathbb{K}$ is a field. If $A$
is an $n\times n$-matrix over $\mathbb{K}$, then the determinant of $A$ can be
computed using (\ref{eq.det.eq.1})... but in practice, you probably do not
\textbf{want} to compute it this way, since the right hand side of
(\ref{eq.det.eq.1}) contains a sum of $n!$ terms.

It turns out that there is an algorithm to compute $\det A$, which is
(usually) a lot faster. It is a version of the Gaussian elimination algorithm
commonly used for solving systems of linear equations.

Let us illustrate it on an example: Set%
\[
n=4,\ \ \ \ \ \ \ \ \ \ \mathbb{K}=\mathbb{Q}\ \ \ \ \ \ \ \ \ \ \text{and
}A=\left(
\begin{array}
[c]{cccc}%
1 & 2 & 3 & 0\\
0 & -1 & 0 & 2\\
2 & 4 & -2 & 3\\
5 & 1 & 3 & 5
\end{array}
\right)  .
\]
We want to find $\det A$.

Exercise \ref{exe.ps4.6k} \textbf{(b)} shows that if we add a scalar multiple
of a column of a matrix to another column of this matrix, then the determinant
of the matrix does not change. Now, by adding appropriate scalar multiples of
the fourth column of $A$ to the first three columns of $A$, we can make sure
that the first three entries of the fourth row of $A$ become zero: Namely, we
have to

\begin{itemize}
\item add $\left(  -1\right)  $ times the fourth column of $A$ to the first
column of $A$;

\item add $\left(  -1/5\right)  $ times the fourth column of $A$ to the second
column of $A$;

\item add $\left(  -3/5\right)  $ times the fourth column of $A$ to the third
column of $A$.
\end{itemize}

These additions can be performed in any order, since none of them
\textquotedblleft interacts\textquotedblright\ with any other (more precisely,
none of them uses any entries that another of them changes). As we know, none
of these additions changes the determinant of the matrix.

Having performed these three additions, we end up with the matrix%
\begin{equation}
A^{\prime}=\left(
\begin{array}
[c]{cccc}%
1 & 2 & 3 & 0\\
-2 & -7/5 & -6/5 & 2\\
-1 & 17/5 & -19/5 & 3\\
0 & 0 & 0 & 5
\end{array}
\right)  . \label{eq.rmk.laplace.pre.gauss-alg.3}%
\end{equation}
We have $\det\left(  A^{\prime}\right)  =\det A$ (because $A^{\prime}$ was
obtained from $A$ by three operations which do not change the determinant).
Moreover, the fourth row of $A^{\prime}$ contains only one nonzero entry --
namely, its last entry. In other words, if we write $A^{\prime}$ in the form
$A^{\prime}=\left(  a_{i,j}^{\prime}\right)  _{1\leq i\leq4,\ 1\leq j\leq4}$,
then $a_{4,j}^{\prime}=0$ for every $j\in\left\{  1,2,3\right\}  $. Thus,
Theorem \ref{thm.laplace.pre} (applied to $4$, $A^{\prime}$ and $a_{i,j}%
^{\prime}$ instead of $n$, $A$ and $a_{i,j}$) shows that%
\begin{align*}
\det\left(  A^{\prime}\right)   &  =\underbrace{a_{4,4}^{\prime}}_{=5}%
\cdot\det\left(  \underbrace{\left(  a_{i,j}^{\prime}\right)  _{1\leq
i\leq3,\ 1\leq j\leq3}}_{=\left(
\begin{array}
[c]{ccc}%
1 & 2 & 3\\
-2 & -7/5 & -6/5\\
-1 & 17/5 & -19/5
\end{array}
\right)  }\right) \\
&  =5\cdot\det\left(
\begin{array}
[c]{ccc}%
1 & 2 & 3\\
-2 & -7/5 & -6/5\\
-1 & 17/5 & -19/5
\end{array}
\right)  .
\end{align*}
Comparing this with $\det\left(  A^{\prime}\right)  =\det A$, we obtain%
\[
\det A=5\cdot\det\left(
\begin{array}
[c]{ccc}%
1 & 2 & 3\\
-2 & -7/5 & -6/5\\
-1 & 17/5 & -19/5
\end{array}
\right)  .
\]

Thus, we have reduced the problem of computing $\det A$ (the determinant of a
$4\times4$-matrix) to the problem of computing $\det\left(
\begin{array}
[c]{ccc}%
1 & 2 & 3\\
-2 & -7/5 & -6/5\\
-1 & 17/5 & -19/5
\end{array}
\right)  $ (the determinant of a $3\times3$-matrix). Likewise, we can try to
reduce the latter problem to the computation of the determinant of a
$2\times2$-matrix, and then further to the computation of the determinant of a
$1\times1$-matrix. (In our example, we obtain $\det A=-140$ at the end.)

This looks like a viable algorithm (which is, furthermore, fairly fast:
essentially as fast as Gaussian elimination). But does it always work? It
turns out that it \textbf{almost} always works. There are cases in which it
can get \textquotedblleft stuck\textquotedblright, and it needs to be modified
to deal with these cases.

Namely, what can happen is that the $\left(  n,n\right)  $-th entry of the
matrix $A$ could be $0$. Again, let us observe this on an example: Set $n=4$
and $A=\left(
\begin{array}
[c]{cccc}%
1 & 2 & 3 & 0\\
0 & -1 & 0 & 2\\
2 & 4 & -2 & 3\\
5 & 1 & 3 & 0
\end{array}
\right)  $. Then, we cannot turn the first three entries of the fourth row of
$A$ into zeroes by adding appropriate multiples of the fourth column to the
first three columns. (Whatever multiples we add, the fourth row stays
unchanged.) However, we can now swap the second row of $A$ with the fourth
row. This operation produces the matrix $B=\left(
\begin{array}
[c]{cccc}%
1 & 2 & 3 & 0\\
5 & 1 & 3 & 0\\
2 & 4 & -2 & 3\\
0 & -1 & 0 & 2
\end{array}
\right)  $, which satisfies $\det B=-\det A$ (by Exercise \ref{exe.ps4.6}
\textbf{(a)}). Thus, it suffices to compute $\det B$; and this can be done as above.

The reason why we swapped the second row of $A$ with the fourth row is that
the last entry of the second row of $A$ was nonzero. In general, we need to
find a $k\in\left\{  1,2,\ldots,n\right\}  $ such that the last entry of the
$k$-th row of $A$ is nonzero, and swap the $k$-th row of $A$ with the $n$-th
row. But what if no such $k$ exists? In this case, we need another way to
compute $\det A$. It turns out that this is very easy: If there is no
$k\in\left\{  1,2,\ldots,n\right\}  $ such that the last entry of the $k$-th
row of $A$ is nonzero, then the last column of $A$ consists of zeroes, and
thus Exercise \ref{exe.ps4.6} \textbf{(d)} shows that $\det A=0$.

When $\mathbb{K}$ is not a field, this algorithm breaks (or, at least,
\textbf{can} break). Indeed, it relies on the fact that the $\left(
n,n\right)  $-th entry of the matrix $A$ is either zero or invertible. Over a
commutative ring $\mathbb{K}$, it might be neither. For example, if we had
tried to work with $\mathbb{K}=\mathbb{Z}$ (instead of $\mathbb{K}=\mathbb{Q}%
$) in our above example, then we would not be able to add $\left(
-1/5\right)  $ times the fourth column of $A$ to the second column of $A$
(because $-1/5\notin\mathbb{Z}=\mathbb{K}$). Fortunately, of course,
$\mathbb{Z}$ is a subset of $\mathbb{Q}$ (and its operations $+$ and $\cdot$
are consistent with those of $\mathbb{Q}$), so that we can just perform the
whole algorithm over $\mathbb{Q}$ instead of $\mathbb{Z}$. However, we aren't
always in luck: Some commutative rings $\mathbb{K}$ cannot be
\textquotedblleft embedded\textquotedblright\ into fields in the way
$\mathbb{Z}$ is embedded into $\mathbb{Q}$. (For instance, $\mathbb{Z}%
/4\mathbb{Z}$ cannot be embedded into a field.)

Nevertheless, there \textbf{are} reasonably fast algorithms for computing
determinants over any commutative ring; see \cite[\S 2]{Rote}.
\end{remark}

\subsection{The Cauchy determinant}

Now, we can state another classical formula for a determinant: the
\textit{Cauchy determinant}. In one of its many forms, it says the following:

\begin{exercise}
\label{exe.cauchy-det}Let $n\in\mathbb{N}$. Let $x_{1},x_{2},\ldots,x_{n}$ be
$n$ elements of $\mathbb{K}$. Let $y_{1},y_{2},\ldots,y_{n}$ be $n$ elements
of $\mathbb{K}$. Assume that $x_{i}+y_{j}$ is invertible in $\mathbb{K}$ for
every $\left(  i,j\right)  \in\left\{  1,2,\ldots,n\right\}  ^{2}$. Then,
prove that%
\[
\det\left(  \left(  \dfrac{1}{x_{i}+y_{j}}\right)  _{1\leq i\leq n,\ 1\leq
j\leq n}\right)  =\frac{\prod_{1\leq i<j\leq n}\left(  \left(  x_{i}%
-x_{j}\right)  \left(  y_{i}-y_{j}\right)  \right)  }{\prod_{\left(
i,j\right)  \in\left\{  1,2,\ldots,n\right\}  ^{2}}\left(  x_{i}+y_{j}\right)
}.
\]

\end{exercise}

There is a different version of the Cauchy determinant floating around in
literature; it differs from Exercise \ref{exe.cauchy-det} in that each
\textquotedblleft$x_{i}+y_{j}$\textquotedblright\ is replaced by
\textquotedblleft$x_{i}-y_{j}$\textquotedblright, and in that
\textquotedblleft$y_{i}-y_{j}$\textquotedblright\ is replaced by
\textquotedblleft$y_{j}-y_{i}$\textquotedblright. Of course, this version is
nothing else than the result of applying Exercise \ref{exe.cauchy-det} to
$-y_{1},-y_{2},\ldots,-y_{n}$ instead of $y_{1},y_{2},\ldots,y_{n}$.

\begin{exercise}
\label{exe.cauchy-det-lem}Let $n$ be a positive integer. Let $\left(
a_{i,j}\right)  _{1\leq i\leq n,\ 1\leq j\leq n}$ be an $n\times n$-matrix
such that $a_{n,n}$ is invertible (in $\mathbb{K}$). Prove that
\begin{align}
&  \det\left(  \left(  a_{i,j}a_{n,n}-a_{i,n}a_{n,j}\right)  _{1\leq i\leq
n-1,\ 1\leq j\leq n-1}\right) \nonumber\\
&  =a_{n,n}^{n-2}\cdot\det\left(  \left(  a_{i,j}\right)  _{1\leq i\leq
n,\ 1\leq j\leq n}\right)  . \label{eq.exe.cauchy-det-lem}%
\end{align}

\end{exercise}

Exercise \ref{exe.cauchy-det-lem} is known as the \textit{Chio pivotal
condensation theorem}\footnote{See \cite[footnote 2]{Heinig} and
\cite[\S 2]{Abeles} for some hints about its history. A variant of the formula
(singling out the $1$-st row and the $1$-st column instead of the $n$-th row
and the $n$-th column) appears in \cite[Chapter Four, Topic \textquotedblleft
Chio's Method\textquotedblright]{Hefferon}.}.

\begin{remark}
Exercise \ref{exe.cauchy-det-lem} gives a way to reduce the computation of an
$n\times n$-determinant (the one on the right hand side of
(\ref{eq.exe.cauchy-det-lem})) to the computation of an $\left(  n-1\right)
\times\left(  n-1\right)  $-determinant (the one on the left hand side),
provided that $a_{n,n}$ is invertible. If this reminds you of Remark
\ref{rmk.laplace.pre.gauss-alg}, you are thinking right...
\end{remark}

\begin{remark}
Exercise \ref{exe.cauchy-det-lem} holds even without the assumption that
$a_{n,n}$ be invertible, as long as we assume (instead) that $n\geq2$. (If we
don't assume that $n\geq2$, then the $a_{n,n}^{n-2}$ on the right hand side of
(\ref{eq.exe.cauchy-det-lem}) will not be defined for non-invertible $a_{n,n}%
$.) Proving this is beyond these notes, though. (A proof of this generalized
version of Exercise \ref{exe.cauchy-det-lem} can be found in \cite{KarZha16}.
It can also be obtained as a particular case of \cite[(4)]{BerBru08}\footnotemark.)
\end{remark}

\footnotetext{In more detail: If we apply \cite[(4)]{BerBru08} to $k=n-1$,
then the right hand side is precisely $\det\left(  \left(  a_{i,j}%
a_{n,n}-a_{i,n}a_{n,j}\right)  _{1\leq i\leq n-1,\ 1\leq j\leq n-1}\right)  $,
and so the formula becomes (\ref{eq.exe.cauchy-det-lem}).}

\subsection{Further determinant equalities}

Next, let us provide an assortment of other exercises on determinants.
Hundreds of exercises (ranging from easy to challenging) on the properties and
evaluations of determinants can be found in \cite[Chapter 2]{FadSom72}, and
some more in \cite[Chapter I]{Prasolov}; in comparison, our selection is
rather small.

\begin{exercise}
\label{exe.det.schur-lem}Let $n\in\mathbb{N}$. Let $\left(  a_{i,j}\right)
_{1\leq i\leq n,\ 1\leq j\leq n}$ be an $n\times n$-matrix. Let $b_{1}$,
$b_{2}$, $\ldots$, $b_{n}$ be $n$ elements of $\mathbb{K}$. Prove that%
\[
\sum\limits_{k=1}^{n}\det\left(  \left(  a_{i,j}b_{i}^{\delta_{j,k}}\right)
_{1\leq i\leq n,\ 1\leq j\leq n}\right)  =\left(  b_{1}+b_{2}+\cdots
+b_{n}\right)  \det\left(  \left(  a_{i,j}\right)  _{1\leq i\leq n,\ 1\leq
j\leq n}\right)  ,
\]
where $\delta_{j,k}$ means the nonnegative integer $%
\begin{cases}
1, & \text{if }j=k;\\
0, & \text{if }j\neq k
\end{cases}
$. Equivalently (in more reader-friendly terms): Prove that%
\begin{align*}
&  \det\left(
\begin{array}
[c]{cccc}%
a_{1,1}b_{1} & a_{1,2} & \cdots & a_{1,n}\\
a_{2,1}b_{2} & a_{2,2} & \cdots & a_{2,n}\\
\vdots & \vdots & \ddots & \vdots\\
a_{n,1}b_{n} & a_{n,2} & \cdots & a_{n,n}%
\end{array}
\right)  +\det\left(
\begin{array}
[c]{cccc}%
a_{1,1} & a_{1,2}b_{1} & \cdots & a_{1,n}\\
a_{2,1} & a_{2,2}b_{2} & \cdots & a_{2,n}\\
\vdots & \vdots & \ddots & \vdots\\
a_{n,1} & a_{n,2}b_{n} & \cdots & a_{n,n}%
\end{array}
\right) \\
&  \ \ \ \ \ \ \ \ \ \ +\cdots+\det\left(
\begin{array}
[c]{cccc}%
a_{1,1} & a_{1,2} & \cdots & a_{1,n}b_{1}\\
a_{2,1} & a_{2,2} & \cdots & a_{2,n}b_{2}\\
\vdots & \vdots & \ddots & \vdots\\
a_{n,1} & a_{n,2} & \cdots & a_{n,n}b_{n}%
\end{array}
\right) \\
&  =\left(  b_{1}+b_{2}+\cdots+b_{n}\right)  \det\left(
\begin{array}
[c]{cccc}%
a_{1,1} & a_{1,2} & \cdots & a_{1,n}\\
a_{2,1} & a_{2,2} & \cdots & a_{2,n}\\
\vdots & \vdots & \ddots & \vdots\\
a_{n,1} & a_{n,2} & \cdots & a_{n,n}%
\end{array}
\right)  .
\end{align*}

\end{exercise}

\begin{exercise}
\label{exe.det.a1a2anx}Let $n\in\mathbb{N}$. Let $a_{1},a_{2},\ldots,a_{n}$ be
$n$ elements of $\mathbb{K}$. Let $x\in\mathbb{K}$. Prove that%
\[
\det\left(
\begin{array}
[c]{cccccc}%
x & a_{1} & a_{2} & \cdots & a_{n-1} & a_{n}\\
a_{1} & x & a_{2} & \cdots & a_{n-1} & a_{n}\\
a_{1} & a_{2} & x & \cdots & a_{n-1} & a_{n}\\
\vdots & \vdots & \vdots & \ddots & \vdots & \vdots\\
a_{1} & a_{2} & a_{3} & \cdots & x & a_{n}\\
a_{1} & a_{2} & a_{3} & \cdots & a_{n} & x
\end{array}
\right)  =\left(  x+\sum_{i=1}^{n}a_{i}\right)  \prod_{i=1}^{n}\left(
x-a_{i}\right)  .
\]

\end{exercise}

\begin{exercise}
\label{exe.det.2diags}Let $n>1$ be an integer. Let $a_{1},a_{2},\ldots,a_{n}$
be $n$ elements of $\mathbb{K}$. Let $b_{1},b_{2},\ldots,b_{n}$ be $n$
elements of $\mathbb{K}$. Let $A$ be the $n\times n$-matrix%
\begin{align*}
&  \left(
\begin{cases}
a_{j}, & \text{if }i=j;\\
b_{j}, & \text{if }i\equiv j+1\operatorname{mod}n;\\
0, & \text{otherwise}%
\end{cases}
\right)  _{1\leq i\leq n,\ 1\leq j\leq n}\\
&  =\left(
\begin{array}
[c]{cccccc}%
a_{1} & 0 & 0 & \cdots & 0 & b_{n}\\
b_{1} & a_{2} & 0 & \cdots & 0 & 0\\
0 & b_{2} & a_{3} & \cdots & 0 & 0\\
\vdots & \vdots & \vdots & \ddots & \vdots & \vdots\\
0 & 0 & 0 & \cdots & a_{n-1} & 0\\
0 & 0 & 0 & \cdots & b_{n-1} & a_{n}%
\end{array}
\right)  .
\end{align*}
Prove that%
\[
\det A=a_{1}a_{2}\cdots a_{n}+\left(  -1\right)  ^{n-1}b_{1}b_{2}\cdots
b_{n}.
\]

\end{exercise}

\begin{remark}
If we replace \textquotedblleft$i\equiv j+1\operatorname{mod}n$%
\textquotedblright\ by \textquotedblleft$i\equiv j+2\operatorname{mod}%
n$\textquotedblright\ in Exercise \ref{exe.det.2diags}, then the pattern can
break. For instance, for $n=4$ we have%
\[
\det\left(
\begin{array}
[c]{cccc}%
a_{1} & 0 & b_{3} & 0\\
0 & a_{2} & 0 & b_{4}\\
b_{1} & 0 & a_{3} & 0\\
0 & b_{2} & 0 & a_{4}%
\end{array}
\right)  =\left(  a_{2}a_{4}-b_{2}b_{4}\right)  \left(  a_{1}a_{3}-b_{1}%
b_{3}\right)  ,
\]
which is not of the form $a_{1}a_{2}a_{3}a_{4}\pm b_{1}b_{2}b_{3}b_{4}$
anymore. Can you guess for which $d\in\left\{  1,2,\ldots,n-1\right\}  $ we
can replace \textquotedblleft$i\equiv j+1\operatorname{mod}n$%
\textquotedblright\ by \textquotedblleft$i\equiv j+d\operatorname{mod}%
n$\textquotedblright\ in Exercise \ref{exe.det.2diags} and still get a formula
of the form $\det A=a_{1}a_{2}\cdots a_{n}\pm b_{1}b_{2}\cdots b_{n}$ ? (The
answer to this question requires a little bit of elementary number theory --
namely, the concept of \textquotedblleft coprimality\textquotedblright. See
\cite[Proposition 6]{HanKra00} for the answer.)
\end{remark}

\subsection{Alternating matrices}

Our next two exercises will concern two special classes of matrices: the
\textit{antisymmetric} and the \textit{alternating matrices}. Let us first
define these classes:

\begin{definition}
\label{def.altern}Let $n\in\mathbb{N}$. Let $A=\left(  a_{i,j}\right)  _{1\leq
i\leq n,\ 1\leq j\leq n}$ be an $n\times n$-matrix.

\textbf{(a)} The matrix $A$ is said to be \textit{antisymmetric} if and only
if $A^{T}=-A$. (Recall that $A^{T}$ is defined as in Definition
\ref{def.transpose}.)

\textbf{(b)} The matrix $A$ is said to be \textit{alternating} if and only if
it satisfies $A^{T}=-A$ and $\left(  a_{i,i}=0\text{ for all }i\in\left\{
1,2,\ldots,n\right\}  \right)  $.
\end{definition}

\begin{example}
\label{exam.altern}A $1\times1$-matrix is alternating if and only if it is the
zero matrix $0_{1\times1}=\left(
\begin{array}
[c]{c}%
0
\end{array}
\right)  $.

A $2\times2$-matrix is alternating if and only if it has the form $\left(
\begin{array}
[c]{cc}%
0 & a\\
-a & 0
\end{array}
\right)  $ for some $a\in\mathbb{K}$.

A $3\times3$-matrix is alternating if and only if it has the form $\left(
\begin{array}
[c]{ccc}%
0 & a & b\\
-a & 0 & c\\
-b & -c & 0
\end{array}
\right)  $ for some $a,b,c\in\mathbb{K}$.

Visually speaking, an $n\times n$-matrix is alternating if and only if its
diagonal entries are $0$ and its entries below the diagonal are the negatives
of their \textquotedblleft mirror-image\textquotedblright\ entries above the diagonal.
\end{example}

\begin{remark}
Clearly, any alternating matrix is antisymmetric. It is easy to see that an
$n\times n$-matrix $A=\left(  a_{i,j}\right)  _{1\leq i\leq n,\ 1\leq j\leq
n}$ is antisymmetric if and only if every $\left(  i,j\right)  \in\left\{
1,2,\ldots,n\right\}  ^{2}$ satisfies $a_{i,j}=-a_{j,i}$. Thus, if $A=\left(
a_{i,j}\right)  _{1\leq i\leq n,\ 1\leq j\leq n}$ is antisymmetric, then every
$i\in\left\{  1,2,\ldots,n\right\}  $ satisfies $a_{i,i}=-a_{i,i}$ and thus
$2a_{i,i}=0$. If $\mathbb{K}$ is one of the rings $\mathbb{Z}$, $\mathbb{Q}$,
$\mathbb{R}$ and $\mathbb{C}$, then we can cancel $2$ from this last equality,
and conclude that every antisymmetric $n\times n$-matrix $A$ is alternating.
However, there are commutative rings $\mathbb{K}$ for which this does not hold
(for example, the ring $\mathbb{Z}/2\mathbb{Z}$ of integers modulo $2$).

Antisymmetric matrices are also known as \textit{skew-symmetric} matrices.
\end{remark}

\begin{exercise}
\label{exe.altern.STAS}Let $n\in\mathbb{N}$ and $m\in\mathbb{N}$. Let $A$ be
an alternating $n\times n$-matrix. Let $S$ be an $n\times m$-matrix. Prove
that the $m\times m$-matrix $S^{T}AS$ is alternating.
\end{exercise}

\begin{exercise}
\label{exe.altern.det}Let $n\in\mathbb{N}$ be odd. Let $A$ be an $n\times
n$-matrix. Prove the following:

\textbf{(a)} If $A$ is antisymmetric, then $2\det A=0$.

\textbf{(b)} If $A$ is alternating, then $\det A=0$.
\end{exercise}

\begin{remark}
If $\mathbb{K}$ is one of the rings $\mathbb{Z}$, $\mathbb{Q}$, $\mathbb{R}$
and $\mathbb{C}$, then Exercise \ref{exe.altern.det} \textbf{(b)} follows from
Exercise \ref{exe.altern.det} \textbf{(a)} (because any alternating matrix is
antisymmetric, and because we can cancel $2$ from the equality $2\det A=0$).
However, this quick way of solving Exercise \ref{exe.altern.det} \textbf{(b)}
does not work for general $\mathbb{K}$.
\end{remark}

\begin{remark}
Exercise \ref{exe.altern.det} \textbf{(b)} provides a really simple formula
for $\det A$ when $A$ is an alternating $n\times n$-matrix for \textbf{odd}
$n$. One might wonder what can be said about $\det A$ when $A$ is an
alternating $n\times n$-matrix for \textbf{even} $n$. The answer is far less
simple, but more interesting: It turns that $\det A$ is the square of a
certain element of $\mathbb{K}$, called the \textit{Pfaffian} of $A$. See
\cite[(5.5)]{Conrad-Pf} for a short introduction into the Pfaffian (although
at a less elementary level than these notes); see \cite[\S 9.5]{BruRys91} and
\cite[\S 12.12]{Loehr-BC}\footnotemark\ for a more combinatorial treatment of
the Pfaffian (and an application to matchings of graphs!). For example, the
Pfaffian of an alternating $4\times4$-matrix $A=\left(
\begin{array}
[c]{cccc}%
0 & a & b & c\\
-a & 0 & d & e\\
-b & -d & 0 & f\\
-c & -e & -f & 0
\end{array}
\right)  $ is $af-be+cd$, and it is indeed easy to check that this matrix
satisfies $\det A=\left(  af-be+cd\right)  ^{2}$.
\end{remark}

\footnotetext{Beware that Loehr, in \cite[\S 12.12]{Loehr-BC}, seems to work
only in the setting where $2$ is cancellable in the ring $\mathbb{K}$ (that
is, where $2a=0$ for an element $a\in\mathbb{K}$ implies $a=0$). Thus, Loehr
does not have to distinguish between antisymmetric and alternating matrices
(he calls them \textquotedblleft skew-symmetric matrices\textquotedblright%
\ instead). His arguments, however, can easily be adapted to the general
case.}

\subsection{Laplace expansion}

We shall now state Laplace expansion in full. We begin with an example:

\begin{example}
\label{exa.laplace.3x3}Let $A=\left(  a_{i,j}\right)  _{1\leq i\leq3,\ 1\leq
j\leq3}$ be a $3\times3$-matrix. From (\ref{eq.det.small.3x3}), we obtain%
\begin{equation}
\det A=a_{1,1}a_{2,2}a_{3,3}+a_{1,2}a_{2,3}a_{3,1}+a_{1,3}a_{2,1}%
a_{3,2}-a_{1,1}a_{2,3}a_{3,2}-a_{1,2}a_{2,1}a_{3,3}-a_{1,3}a_{2,2}a_{3,1}.
\label{eq.exa.laplace.3x3.1}%
\end{equation}
On the right hand side of this equality, we have six terms, each of which
contains either $a_{2,1}$ or $a_{2,2}$ or $a_{2,3}$. Let us combine the two
terms containing $a_{2,1}$ and factor out $a_{2,1}$, then do the same with the
two terms containing $a_{2,2}$, and with the two terms containing $a_{2,3}$.
As a result, (\ref{eq.exa.laplace.3x3.1}) becomes%
\begin{align}
&  \det A\nonumber\\
&  =a_{1,1}a_{2,2}a_{3,3}+a_{1,2}a_{2,3}a_{3,1}+a_{1,3}a_{2,1}a_{3,2}%
-a_{1,1}a_{2,3}a_{3,2}-a_{1,2}a_{2,1}a_{3,3}-a_{1,3}a_{2,2}a_{3,1}\nonumber\\
&  =a_{2,1}\underbrace{\left(  a_{1,3}a_{3,2}-a_{1,2}a_{3,3}\right)  }%
_{=\det\left(
\begin{array}
[c]{cc}%
a_{1,3} & a_{1,2}\\
a_{3,3} & a_{3,2}%
\end{array}
\right)  }+a_{2,2}\underbrace{\left(  a_{1,1}a_{3,3}-a_{1,3}a_{3,1}\right)
}_{=\det\left(
\begin{array}
[c]{cc}%
a_{1,1} & a_{1,3}\\
a_{3,1} & a_{3,3}%
\end{array}
\right)  }+a_{2,3}\underbrace{\left(  a_{1,2}a_{3,1}-a_{1,1}a_{3,2}\right)
}_{=\det\left(
\begin{array}
[c]{cc}%
a_{1,2} & a_{1,1}\\
a_{3,2} & a_{3,1}%
\end{array}
\right)  }\nonumber\\
&  =a_{2,1}\det\left(
\begin{array}
[c]{cc}%
a_{1,3} & a_{1,2}\\
a_{3,3} & a_{3,2}%
\end{array}
\right)  +a_{2,2}\det\left(
\begin{array}
[c]{cc}%
a_{1,1} & a_{1,3}\\
a_{3,1} & a_{3,3}%
\end{array}
\right)  +a_{2,3}\det\left(
\begin{array}
[c]{cc}%
a_{1,2} & a_{1,1}\\
a_{3,2} & a_{3,1}%
\end{array}
\right)  . \label{eq.exa.laplace.3x3.2}%
\end{align}
This is a nice formula with an obvious pattern: The right hand side can be
rewritten as $\sum_{q=1}^{3}a_{2,q}\det\left(  B_{2,q}\right)  $, where
$B_{2,q}=\left(
\begin{array}
[c]{cc}%
a_{1,q+2} & a_{1,q+1}\\
a_{3,q+2} & a_{3,q+1}%
\end{array}
\right)  $ (where we set $a_{i,4}=a_{i,1}$ and $a_{i,5}=a_{i,2}$ for all
$i\in\left\{  1,2,3\right\}  $). Notice the cyclic symmetry (with respect to
the index of the column) in this formula! Unfortunately, in this exact form,
the formula does not generalize to bigger matrices (or even to smaller: the
analogue for a $2\times2$-matrix would be $\det\left(
\begin{array}
[c]{cc}%
a_{1,1} & a_{1,2}\\
a_{2,1} & a_{2,2}%
\end{array}
\right)  =-a_{2,1}a_{1,2}+a_{2,2}a_{1,1}$, which has a minus sign unlike
$\sum_{q=1}^{3}a_{2,q}\det\left(  B_{2,q}\right)  $).

However, we can slightly modify our formula, sacrificing the cyclic symmetry
but making it generalize. Namely, let us rewrite $a_{1,3}a_{3,2}%
-a_{1,2}a_{3,3}$ as $-\left(  a_{1,2}a_{3,3}-a_{1,3}a_{3,2}\right)  $ and
$a_{1,2}a_{3,1}-a_{1,1}a_{3,2}$ as $-\left(  a_{1,1}a_{3,2}-a_{1,2}%
a_{3,1}\right)  $; we thus obtain
\begin{align}
&  \det A\nonumber\\
&  =a_{2,1}\underbrace{\left(  a_{1,3}a_{3,2}-a_{1,2}a_{3,3}\right)
}_{=-\left(  a_{1,2}a_{3,3}-a_{1,3}a_{3,2}\right)  }+a_{2,2}\left(
a_{1,1}a_{3,3}-a_{1,3}a_{3,1}\right)  +a_{2,3}\underbrace{\left(
a_{1,2}a_{3,1}-a_{1,1}a_{3,2}\right)  }_{=-\left(  a_{1,1}a_{3,2}%
-a_{1,2}a_{3,1}\right)  }\nonumber\\
&  =-a_{2,1}\underbrace{\left(  a_{1,2}a_{3,3}-a_{1,3}a_{3,2}\right)  }%
_{=\det\left(
\begin{array}
[c]{cc}%
a_{1,2} & a_{1,3}\\
a_{3,2} & a_{3,3}%
\end{array}
\right)  }+a_{2,2}\underbrace{\left(  a_{1,1}a_{3,3}-a_{1,3}a_{3,1}\right)
}_{=\det\left(
\begin{array}
[c]{cc}%
a_{1,1} & a_{1,3}\\
a_{3,1} & a_{3,3}%
\end{array}
\right)  }-a_{2,3}\underbrace{\left(  a_{1,1}a_{3,2}-a_{1,2}a_{3,1}\right)
}_{=\det\left(
\begin{array}
[c]{cc}%
a_{1,1} & a_{1,2}\\
a_{3,1} & a_{3,2}%
\end{array}
\right)  }\nonumber\\
&  =-a_{2,1}\det\left(
\begin{array}
[c]{cc}%
a_{1,2} & a_{1,3}\\
a_{3,2} & a_{3,3}%
\end{array}
\right)  +a_{2,2}\det\left(
\begin{array}
[c]{cc}%
a_{1,1} & a_{1,3}\\
a_{3,1} & a_{3,3}%
\end{array}
\right)  -a_{2,3}\det\left(
\begin{array}
[c]{cc}%
a_{1,1} & a_{1,2}\\
a_{3,1} & a_{3,2}%
\end{array}
\right) \nonumber\\
&  =\sum_{q=1}^{3}\left(  -1\right)  ^{q}a_{2,q}\det\left(  C_{2,q}\right)  ,
\label{eq.exa.laplace.3x3.3}%
\end{align}
where $C_{2,q}$ means the matrix obtained from $A$ by crossing out the $2$-nd
row and the $q$-th column. This formula (unlike (\ref{eq.exa.laplace.3x3.2}))
involves powers of $-1$, but it can be generalized.

How? First, we notice that we can find a similar formula by factoring out
$a_{1,1},a_{1,2},a_{1,3}$ (instead of $a_{2,1},a_{2,2},a_{2,3}$); this formula
will be%
\[
\det A=\sum_{q=1}^{3}\left(  -1\right)  ^{q-1}a_{1,q}\det\left(
C_{1,q}\right)  ,
\]
where $C_{1,q}$ means the matrix obtained from $A$ by crossing out the $1$-st
row and the $q$-th column. This formula, and (\ref{eq.exa.laplace.3x3.3}),
suggest the following generalization: If $A=\left(  a_{i,j}\right)  _{1\leq
i\leq n,\ 1\leq j\leq n}$ is an $n\times n$-matrix, and if $p\in\left\{
1,2,\ldots,n\right\}  $, then%
\begin{equation}
\det A=\sum_{q=1}^{n}\left(  -1\right)  ^{p+q}a_{p,q}\det\left(
C_{p,q}\right)  , \label{eq.exa.laplace.3x3.4}%
\end{equation}
where $C_{p,q}$ means the matrix obtained from $A$ by crossing out the $p$-th
row and the $q$-th column. (The only part of this formula which is not easy to
guess is $\left(  -1\right)  ^{p+q}$; you might need to compute several
particular cases to guess this pattern. Of course, you could also have guessed
$\left(  -1\right)  ^{p-q}$ or $\left(  -1\right)  ^{q-p}$ instead, because
$\left(  -1\right)  ^{p+q}=\left(  -1\right)  ^{p-q}=\left(  -1\right)
^{q-p}$.)

The formula (\ref{eq.exa.laplace.3x3.4}) is what is usually called the Laplace
expansion with respect to the $p$-th row. We will prove it below (Theorem
\ref{thm.laplace.gen} \textbf{(a)}), and we will also prove an analogous
\textquotedblleft Laplace expansion with respect to the $q$-th
column\textquotedblright\ (Theorem \ref{thm.laplace.gen} \textbf{(b)}).
\end{example}

Let us first define a notation:

\begin{definition}
\label{def.submatrix}Let $n\in\mathbb{N}$ and $m\in\mathbb{N}$. Let $A=\left(
a_{i,j}\right)  _{1\leq i\leq n,\ 1\leq j\leq m}$ be an $n\times m$-matrix.
Let $i_{1},i_{2},\ldots,i_{u}$ be some elements of $\left\{  1,2,\ldots
,n\right\}  $; let $j_{1},j_{2},\ldots,j_{v}$ be some elements of $\left\{
1,2,\ldots,m\right\}  $. Then, we define $\operatorname*{sub}\nolimits_{i_{1}%
,i_{2},\ldots,i_{u}}^{j_{1},j_{2},\ldots,j_{v}}A$ to be the $u\times v$-matrix
$\left(  a_{i_{x},j_{y}}\right)  _{1\leq x\leq u,\ 1\leq y\leq v}$.

When $i_{1}<i_{2}<\cdots<i_{u}$ and $j_{1}<j_{2}<\cdots<j_{v}$, the matrix
$\operatorname*{sub}\nolimits_{i_{1},i_{2},\ldots,i_{u}}^{j_{1},j_{2}%
,\ldots,j_{v}}A$ can be obtained from $A$ by crossing out all rows other than
the $i_{1}$-th, the $i_{2}$-th, etc., the $i_{u}$-th row and crossing out all
columns other than the $j_{1}$-th, the $j_{2}$-th, etc., the $j_{v}$-th
column. Thus, in this case, $\operatorname*{sub}\nolimits_{i_{1},i_{2}%
,\ldots,i_{u}}^{j_{1},j_{2},\ldots,j_{v}}A$ is called a \textit{submatrix} of
$A$.
\end{definition}

For example, if $n=3$, $m=4$ and $A=\left(
\begin{array}
[c]{cccc}%
a & b & c & d\\
e & f & g & h\\
i & j & k & \ell
\end{array}
\right)  $, then $\operatorname*{sub}\nolimits_{1,3}^{2,3,4}A=\left(
\begin{array}
[c]{ccc}%
b & c & d\\
j & k & \ell
\end{array}
\right)  $ (this is a submatrix of $A$) and $\operatorname*{sub}%
\nolimits_{2,3}^{3,1,1}A=\left(
\begin{array}
[c]{ccc}%
g & e & e\\
k & i & i
\end{array}
\right)  $ (this is not, in general, a submatrix of $A$).

The following properties follow trivially from the definitions:

\begin{proposition}
\label{prop.submatrix.easy}Let $n\in\mathbb{N}$ and $m\in\mathbb{N}$. Let $A$
be an $n\times m$-matrix. Recall the notations introduced in Definition
\ref{def.rowscols}.

\textbf{(a)} We have $\operatorname*{sub}\nolimits_{1,2,\ldots,n}%
^{1,2,\ldots,m}A=A$.

\textbf{(b)} If $i_{1},i_{2},\ldots,i_{u}$ are some elements of $\left\{
1,2,\ldots,n\right\}  $, then%
\[
\operatorname*{rows}\nolimits_{i_{1},i_{2},\ldots,i_{u}}A=\operatorname*{sub}%
\nolimits_{i_{1},i_{2},\ldots,i_{u}}^{1,2,\ldots,m}A.
\]

\textbf{(c)} If $j_{1},j_{2},\ldots,j_{v}$ are some elements of $\left\{
1,2,\ldots,m\right\}  $, then%
\[
\operatorname*{cols}\nolimits_{j_{1},j_{2},\ldots,j_{v}}A=\operatorname*{sub}%
\nolimits_{1,2,\ldots,n}^{j_{1},j_{2},\ldots,j_{v}}A.
\]

\textbf{(d)} Let $i_{1},i_{2},\ldots,i_{u}$ be some elements of $\left\{
1,2,\ldots,n\right\}  $; let $j_{1},j_{2},\ldots,j_{v}$ be some elements of
$\left\{  1,2,\ldots,m\right\}  $. Then,%
\[
\operatorname*{sub}\nolimits_{i_{1},i_{2},\ldots,i_{u}}^{j_{1},j_{2}%
,\ldots,j_{v}}A=\operatorname*{rows}\nolimits_{i_{1},i_{2},\ldots,i_{u}%
}\left(  \operatorname*{cols}\nolimits_{j_{1},j_{2},\ldots,j_{v}}A\right)
=\operatorname*{cols}\nolimits_{j_{1},j_{2},\ldots,j_{v}}\left(
\operatorname*{rows}\nolimits_{i_{1},i_{2},\ldots,i_{u}}A\right)  .
\]

\textbf{(e)} Let $i_{1},i_{2},\ldots,i_{u}$ be some elements of $\left\{
1,2,\ldots,n\right\}  $; let $j_{1},j_{2},\ldots,j_{v}$ be some elements of
$\left\{  1,2,\ldots,m\right\}  $. Then,%
\[
\left(  \operatorname*{sub}\nolimits_{i_{1},i_{2},\ldots,i_{u}}^{j_{1}%
,j_{2},\ldots,j_{v}}A\right)  ^{T}=\operatorname*{sub}\nolimits_{j_{1}%
,j_{2},\ldots,j_{v}}^{i_{1},i_{2},\ldots,i_{u}}\left(  A^{T}\right)  .
\]

\end{proposition}

\begin{verlong}
\begin{proof}
[Proof of Proposition \ref{prop.submatrix.easy}.]Write the matrix $A$ in the
form $A=\left(  a_{i,j}\right)  _{1\leq i\leq n,\ 1\leq j\leq m}$.

\textbf{(a)} We have $A=\left(  a_{i,j}\right)  _{1\leq i\leq n,\ 1\leq j\leq
m}$. Thus, the definition of $\operatorname*{sub}\nolimits_{1,2,\ldots
,n}^{1,2,\ldots,m}A$ yields%
\begin{align*}
\operatorname*{sub}\nolimits_{1,2,\ldots,n}^{1,2,\ldots,m}A  &  =\left(
a_{x,y}\right)  _{1\leq x\leq n,\ 1\leq y\leq m}=\left(  a_{i,j}\right)
_{1\leq i\leq n,\ 1\leq j\leq m}\\
&  \ \ \ \ \ \ \ \ \ \ \left(  \text{here, we renamed the index }\left(
x,y\right)  \text{ as }\left(  i,j\right)  \right) \\
&  =A.
\end{align*}
This proves Proposition \ref{prop.submatrix.easy} \textbf{(a)}.

\textbf{(b)} Let $i_{1},i_{2},\ldots,i_{u}$ be some elements of $\left\{
1,2,\ldots,n\right\}  $. The definition of $\operatorname*{rows}%
\nolimits_{i_{1},i_{2},\ldots,i_{u}}A$ yields $\operatorname*{rows}%
\nolimits_{i_{1},i_{2},\ldots,i_{u}}A=\left(  a_{i_{x},j}\right)  _{1\leq
x\leq u,\ 1\leq j\leq m}$ (since $A=\left(  a_{i,j}\right)  _{1\leq i\leq
n,\ 1\leq j\leq m}$). On the other hand, we have $A=\left(  a_{i,j}\right)
_{1\leq i\leq n,\ 1\leq j\leq m}$. Hence, the definition of
$\operatorname*{sub}\nolimits_{i_{1},i_{2},\ldots,i_{u}}^{1,2,\ldots,m}A$
yields
\[
\operatorname*{sub}\nolimits_{i_{1},i_{2},\ldots,i_{u}}^{1,2,\ldots
,m}A=\left(  a_{i_{x},y}\right)  _{1\leq x\leq u,\ 1\leq y\leq m}=\left(
a_{i_{x},j}\right)  _{1\leq x\leq u,\ 1\leq j\leq m}%
\]
(here, we renamed the index $\left(  x,y\right)  $ as $\left(  x,j\right)  $).
Comparing this with $\operatorname*{rows}\nolimits_{i_{1},i_{2},\ldots,i_{u}%
}A=\left(  a_{i_{x},j}\right)  _{1\leq x\leq u,\ 1\leq j\leq m}$, we obtain
$\operatorname*{rows}\nolimits_{i_{1},i_{2},\ldots,i_{u}}A=\operatorname*{sub}%
\nolimits_{i_{1},i_{2},\ldots,i_{u}}^{1,2,\ldots,m}A$. This proves Proposition
\ref{prop.submatrix.easy} \textbf{(b)}.

\textbf{(c)} Let $j_{1},j_{2},\ldots,j_{v}$ be some elements of $\left\{
1,2,\ldots,m\right\}  $. The definition of $\operatorname*{cols}%
\nolimits_{j_{1},j_{2},\ldots,j_{v}}A$ yields $\operatorname*{cols}%
\nolimits_{j_{1},j_{2},\ldots,j_{v}}A=\left(  a_{i,j_{y}}\right)  _{1\leq
i\leq n,\ 1\leq y\leq v}$ (since $A=\left(  a_{i,j}\right)  _{1\leq i\leq
n,\ 1\leq j\leq m}$). On the other hand, we have $A=\left(  a_{i,j}\right)
_{1\leq i\leq n,\ 1\leq j\leq m}$. Hence, the definition of
$\operatorname*{sub}\nolimits_{1,2,\ldots,n}^{j_{1},j_{2},\ldots,j_{v}}A$
yields%
\[
\operatorname*{sub}\nolimits_{1,2,\ldots,n}^{j_{1},j_{2},\ldots,j_{v}%
}A=\left(  a_{x,j_{y}}\right)  _{1\leq x\leq n,\ 1\leq y\leq v}=\left(
a_{i,j_{y}}\right)  _{1\leq i\leq n,\ 1\leq y\leq v}%
\]
(here, we renamed the index $\left(  x,y\right)  $ as $\left(  i,y\right)  $).
Comparing this with $\operatorname*{cols}\nolimits_{j_{1},j_{2},\ldots,j_{v}%
}A=\left(  a_{i,j_{y}}\right)  _{1\leq i\leq n,\ 1\leq y\leq v}$, we obtain
$\operatorname*{cols}\nolimits_{j_{1},j_{2},\ldots,j_{v}}A=\operatorname*{sub}%
\nolimits_{1,2,\ldots,n}^{j_{1},j_{2},\ldots,j_{v}}A$. This proves Proposition
\ref{prop.submatrix.easy} \textbf{(c)}.

\textbf{(d)} We have $A=\left(  a_{i,j}\right)  _{1\leq i\leq n,\ 1\leq j\leq
m}$. Thus, the definition of $\operatorname*{sub}\nolimits_{i_{1},i_{2}%
,\ldots,i_{u}}^{j_{1},j_{2},\ldots,j_{v}}A$ yields $\operatorname*{sub}%
\nolimits_{i_{1},i_{2},\ldots,i_{u}}^{j_{1},j_{2},\ldots,j_{v}}A=\left(
a_{i_{x},j_{y}}\right)  _{1\leq x\leq u,\ 1\leq y\leq v}$.

On the other hand, we have $A=\left(  a_{i,j}\right)  _{1\leq i\leq n,\ 1\leq
j\leq m}$. Thus, the definition of $\operatorname*{cols}\nolimits_{j_{1}%
,j_{2},\ldots,j_{v}}A$ yields
\[
\operatorname*{cols}\nolimits_{j_{1},j_{2},\ldots,j_{v}}A=\left(  a_{i,j_{y}%
}\right)  _{1\leq i\leq n,\ 1\leq y\leq v}=\left(  a_{i,j_{j}}\right)  _{1\leq
i\leq n,\ 1\leq j\leq v}%
\]
(here, we renamed the index $\left(  i,y\right)  $ as $\left(  i,j\right)
$)\ \ \ \ \footnote{The notation $j_{j}$ might look fishy, since it uses the
letter \textquotedblleft$j$\textquotedblright\ in two unrelated meanings: Its
meaning in \textquotedblleft$j_{1},j_{2},\ldots,j_{v}$\textquotedblright\ has
nothing to do with its meaning with \textquotedblleft$1\leq j\leq
v$\textquotedblright. However, it is easy to distinguish between these two
kinds of \textquotedblleft$j$\textquotedblright, because the former always
appears with a subscript, whereas the latter never does.}. Hence, the
definition of \newline$\operatorname*{rows}\nolimits_{i_{1},i_{2},\ldots
,i_{u}}\left(  \operatorname*{cols}\nolimits_{j_{1},j_{2},\ldots,j_{v}%
}A\right)  $ yields%
\[
\operatorname*{rows}\nolimits_{i_{1},i_{2},\ldots,i_{u}}\left(
\operatorname*{cols}\nolimits_{j_{1},j_{2},\ldots,j_{v}}A\right)  =\left(
a_{i_{x},j_{j}}\right)  _{1\leq x\leq n,\ 1\leq j\leq v}=\left(
a_{i_{x},j_{y}}\right)  _{1\leq x\leq n,\ 1\leq y\leq v}%
\]
(here, we renamed the index $\left(  x,j\right)  $ as $\left(  x,y\right)  $).
Compared with $\operatorname*{sub}\nolimits_{i_{1},i_{2},\ldots,i_{u}}%
^{j_{1},j_{2},\ldots,j_{v}}A=\left(  a_{i_{x},j_{y}}\right)  _{1\leq x\leq
u,\ 1\leq y\leq v}$, this yields%
\begin{equation}
\operatorname*{sub}\nolimits_{i_{1},i_{2},\ldots,i_{u}}^{j_{1},j_{2}%
,\ldots,j_{v}}A=\operatorname*{rows}\nolimits_{i_{1},i_{2},\ldots,i_{u}%
}\left(  \operatorname*{cols}\nolimits_{j_{1},j_{2},\ldots,j_{v}}A\right)  .
\label{pf.prop.submatrix.easy.d.1}%
\end{equation}

Furthermore, we have $A=\left(  a_{i,j}\right)  _{1\leq i\leq n,\ 1\leq j\leq
m}$. Thus, the definition of $\operatorname*{rows}\nolimits_{i_{1}%
,i_{2},\ldots,i_{u}}A$ yields
\[
\operatorname*{rows}\nolimits_{i_{1},i_{2},\ldots,i_{u}}A=\left(  a_{i_{x}%
,j}\right)  _{1\leq x\leq u,\ 1\leq j\leq m}=\left(  a_{i_{i},j}\right)
_{1\leq i\leq u,\ 1\leq j\leq m}%
\]
\footnote{The double use of the letter \textquotedblleft$i$\textquotedblright%
\ in \textquotedblleft$i_{i}$\textquotedblright\ might appear confusing. The
first \textquotedblleft$i$\textquotedblright\ is part of the notation $i_{k}$
for $k\in\left\{  1,2,\ldots,u\right\}  $; the second \textquotedblleft%
$i$\textquotedblright\ is an element of $\left\{  1,2,\ldots,u\right\}  $.
These two \textquotedblleft$i$\textquotedblright s are unrelated to each
other. I hope the reader can easily tell them apart by the fact that the
\textquotedblleft$i$\textquotedblright\ that is part of the notation $i_{k}$
always appears with a subscript, whereas the second \textquotedblleft%
$i$\textquotedblright\ never does.} (here, we renamed the index $\left(
x,j\right)  $ as $\left(  i,j\right)  $). Hence, the definition of
\newline$\operatorname*{cols}\nolimits_{j_{1},j_{2},\ldots,j_{v}}\left(
\operatorname*{rows}\nolimits_{i_{1},i_{2},\ldots,i_{u}}A\right)  $ yields%
\[
\operatorname*{cols}\nolimits_{j_{1},j_{2},\ldots,j_{v}}\left(
\operatorname*{rows}\nolimits_{i_{1},i_{2},\ldots,i_{u}}A\right)  =\left(
a_{i_{i},j_{y}}\right)  _{1\leq i\leq n,\ 1\leq y\leq v}=\left(
a_{i_{x},j_{y}}\right)  _{1\leq x\leq n,\ 1\leq y\leq v}%
\]
(here, we renamed the index $\left(  i,y\right)  $ as $\left(  x,y\right)  $).
Compared with $\operatorname*{sub}\nolimits_{i_{1},i_{2},\ldots,i_{u}}%
^{j_{1},j_{2},\ldots,j_{v}}A=\left(  a_{i_{x},j_{y}}\right)  _{1\leq x\leq
u,\ 1\leq y\leq v}$, this yields%
\[
\operatorname*{sub}\nolimits_{i_{1},i_{2},\ldots,i_{u}}^{j_{1},j_{2}%
,\ldots,j_{v}}A=\operatorname*{cols}\nolimits_{j_{1},j_{2},\ldots,j_{v}%
}\left(  \operatorname*{rows}\nolimits_{i_{1},i_{2},\ldots,i_{u}}A\right)  .
\]
Combining this with (\ref{pf.prop.submatrix.easy.d.1}), we obtain%
\[
\operatorname*{sub}\nolimits_{i_{1},i_{2},\ldots,i_{u}}^{j_{1},j_{2}%
,\ldots,j_{v}}A=\operatorname*{rows}\nolimits_{i_{1},i_{2},\ldots,i_{u}%
}\left(  \operatorname*{cols}\nolimits_{j_{1},j_{2},\ldots,j_{v}}A\right)
=\operatorname*{cols}\nolimits_{j_{1},j_{2},\ldots,j_{v}}\left(
\operatorname*{rows}\nolimits_{i_{1},i_{2},\ldots,i_{u}}A\right)  .
\]
This proves Proposition \ref{prop.submatrix.easy} \textbf{(d)}.

\textbf{(e)} We have $A=\left(  a_{i,j}\right)  _{1\leq i\leq n,\ 1\leq j\leq
m}$. Thus, the definition of $\operatorname*{sub}\nolimits_{i_{1},i_{2}%
,\ldots,i_{u}}^{j_{1},j_{2},\ldots,j_{v}}A$ yields
\[
\operatorname*{sub}\nolimits_{i_{1},i_{2},\ldots,i_{u}}^{j_{1},j_{2}%
,\ldots,j_{v}}A=\left(  a_{i_{x},j_{y}}\right)  _{1\leq x\leq u,\ 1\leq y\leq
v}=\left(  a_{i_{i},j_{j}}\right)  _{1\leq i\leq u,\ 1\leq j\leq v}%
\]
(here, we have renamed the index $\left(  x,y\right)  $ as $\left(
i,j\right)  $). Therefore, the definition of $\left(  \operatorname*{sub}%
\nolimits_{i_{1},i_{2},\ldots,i_{u}}^{j_{1},j_{2},\ldots,j_{v}}A\right)  ^{T}$
yields
\[
\left(  \operatorname*{sub}\nolimits_{i_{1},i_{2},\ldots,i_{u}}^{j_{1}%
,j_{2},\ldots,j_{v}}A\right)  ^{T}=\left(  a_{i_{j},j_{i}}\right)  _{1\leq
i\leq v,\ 1\leq j\leq u}=\left(  a_{i_{y},j_{x}}\right)  _{1\leq x\leq
v,\ 1\leq y\leq u}%
\]
(here, we have renamed the index $\left(  i,j\right)  $ as $\left(
x,y\right)  $).

On the other hand, the definition of $A^{T}$ yields $A^{T}=\left(
a_{j,i}\right)  _{1\leq i\leq m,\ 1\leq j\leq n}$ (since $A=\left(
a_{i,j}\right)  _{1\leq i\leq n,\ 1\leq j\leq m}$). Hence, the definition of
$\operatorname*{sub}\nolimits_{j_{1},j_{2},\ldots,j_{v}}^{i_{1},i_{2}%
,\ldots,i_{u}}\left(  A^{T}\right)  $ yields $\operatorname*{sub}%
\nolimits_{j_{1},j_{2},\ldots,j_{v}}^{i_{1},i_{2},\ldots,i_{u}}\left(
A^{T}\right)  =\left(  a_{i_{y},j_{x}}\right)  _{1\leq x\leq v,\ 1\leq y\leq
u}$. Comparing this with $\left(  \operatorname*{sub}\nolimits_{i_{1}%
,i_{2},\ldots,i_{u}}^{j_{1},j_{2},\ldots,j_{v}}A\right)  ^{T}=\left(
a_{i_{y},j_{x}}\right)  _{1\leq x\leq v,\ 1\leq y\leq u}$, we obtain%
\[
\left(  \operatorname*{sub}\nolimits_{i_{1},i_{2},\ldots,i_{u}}^{j_{1}%
,j_{2},\ldots,j_{v}}A\right)  ^{T}=\operatorname*{sub}\nolimits_{j_{1}%
,j_{2},\ldots,j_{v}}^{i_{1},i_{2},\ldots,i_{u}}\left(  A^{T}\right)  .
\]
This proves Proposition \ref{prop.submatrix.easy} \textbf{(e)}.
\end{proof}
\end{verlong}

\begin{definition}
\label{def.hat-omit}Let $n\in\mathbb{N}$. Let $a_{1},a_{2},\ldots,a_{n}$ be
$n$ objects. Let $i\in\left\{  1,2,\ldots,n\right\}  $. Then, $\left(
a_{1},a_{2},\ldots,\widehat{a_{i}},\ldots,a_{n}\right)  $ shall mean the list
$\left(  a_{1},a_{2},\ldots,a_{i-1},a_{i+1},a_{i+2},\ldots,a_{n}\right)  $
(that is, the list $\left(  a_{1},a_{2},\ldots,a_{n}\right)  $ with its $i$-th
entry removed). (Thus, the \textquotedblleft hat\textquotedblright\ over the
$a_{i}$ means that this $a_{i}$ is being omitted from the list.)

For example, $\left(  1^{2},2^{2},\ldots,\widehat{5^{2}},\ldots,8^{2}\right)
=\left(  1^{2},2^{2},3^{2},4^{2},6^{2},7^{2},8^{2}\right)  $.
\end{definition}

\begin{definition}
\label{def.submatrix.minor}Let $n\in\mathbb{N}$ and $m\in\mathbb{N}$. Let $A$
be an $n\times m$-matrix. For every $i\in\left\{  1,2,\ldots,n\right\}  $ and
$j\in\left\{  1,2,\ldots,m\right\}  $, we let $A_{\sim i,\sim j}$ be the
$\left(  n-1\right)  \times\left(  m-1\right)  $-matrix $\operatorname*{sub}%
\nolimits_{1,2,\ldots,\widehat{i},\ldots,n}^{1,2,\ldots,\widehat{j},\ldots
,m}A$. (Thus, $A_{\sim i,\sim j}$ is the matrix obtained from $A$ by crossing
out the $i$-th row and the $j$-th column.)

For example, if $n=m=3$ and $A=\left(
\begin{array}
[c]{ccc}%
a & b & c\\
d & e & f\\
g & h & i
\end{array}
\right)  $, then $A_{\sim1,\sim2}=\left(
\begin{array}
[c]{cc}%
d & f\\
g & i
\end{array}
\right)  $ and $A_{\sim3,\sim2}=\left(
\begin{array}
[c]{cc}%
a & c\\
d & f
\end{array}
\right)  $.
\end{definition}

The notation $A_{\sim i,\sim j}$ introduced in Definition
\ref{def.submatrix.minor} is not very standard; but there does not seem to be
a standard one\footnote{For example, Gill Williamson uses the notation
$A\left(  i\mid j\right)  $ in \cite[Chapter 3]{Gill}.}.

Now we can finally state Laplace expansion:

\begin{theorem}
\label{thm.laplace.gen}Let $n\in\mathbb{N}$. Let $A=\left(  a_{i,j}\right)
_{1\leq i\leq n,\ 1\leq j\leq n}$ be an $n\times n$-matrix.

\textbf{(a)} For every $p\in\left\{  1,2,\ldots,n\right\}  $, we have%
\[
\det A=\sum_{q=1}^{n}\left(  -1\right)  ^{p+q}a_{p,q}\det\left(  A_{\sim
p,\sim q}\right)  .
\]

\textbf{(b)} For every $q\in\left\{  1,2,\ldots,n\right\}  $, we have%
\[
\det A=\sum_{p=1}^{n}\left(  -1\right)  ^{p+q}a_{p,q}\det\left(  A_{\sim
p,\sim q}\right)  .
\]

\end{theorem}

Theorem \ref{thm.laplace.gen} \textbf{(a)} is known as the \textit{Laplace
expansion along the }$p$\textit{-th row} (or \textit{Laplace expansion with
respect to the }$p$\textit{-th row}), whereas Theorem \ref{thm.laplace.gen}
\textbf{(b)} is known as the \textit{Laplace expansion along the }%
$q$\textit{-th column} (or \textit{Laplace expansion with respect to the }%
$q$\textit{-th column}). Notice that Theorem \ref{thm.laplace.gen}
\textbf{(a)} is equivalent to the formula (\ref{eq.exa.laplace.3x3.4}),
because the $A_{\sim p,\sim q}$ in Theorem \ref{thm.laplace.gen} \textbf{(a)}
is precisely what we called $C_{p,q}$ in (\ref{eq.exa.laplace.3x3.4}).

We prepare the field for the proof of Theorem \ref{thm.laplace.gen} with a few lemmas.

\begin{vershort}
\begin{lemma}
\label{lem.laplace.gpshort}For every $n\in\mathbb{N}$, let $\left[  n\right]
$ denote the set $\left\{  1,2,\ldots,n\right\}  $.

Let $n\in\mathbb{N}$. For every $p\in\left[  n\right]  $, we define a
permutation $g_{p}\in S_{n}$ by $g_{p}=\operatorname*{cyc}%
\nolimits_{p,p+1,\ldots,n}$ (where we are using the notations of Definition
\ref{def.perm.cycles}).

\textbf{(a)} We have $\left(  g_{p}\left(  1\right)  ,g_{p}\left(  2\right)
,\ldots,g_{p}\left(  n-1\right)  \right)  =\left(  1,2,\ldots,\widehat{p}%
,\ldots,n\right)  $ for every $p\in\left[  n\right]  $.

\textbf{(b)} We have $\left(  -1\right)  ^{g_{p}}=\left(  -1\right)  ^{n-p}$
for every $p\in\left[  n\right]  $.

\textbf{(c)} Let $p\in\left[  n\right]  $. We define a map%
\[
g_{p}^{\prime}:\left[  n-1\right]  \rightarrow\left[  n\right]  \setminus
\left\{  p\right\}
\]
by%
\[
\left(  g_{p}^{\prime}\left(  i\right)  =g_{p}\left(  i\right)
\ \ \ \ \ \ \ \ \ \ \text{for every }i\in\left[  n-1\right]  \right)  .
\]
This map $g_{p}^{\prime}$ is well-defined and bijective.

\textbf{(d)} Let $p\in\left[  n\right]  $ and $q\in\left[  n\right]  $. We
define a map
\[
T:\left\{  \tau\in S_{n}\ \mid\ \tau\left(  n\right)  =n\right\}
\rightarrow\left\{  \tau\in S_{n}\ \mid\ \tau\left(  p\right)  =q\right\}
\]
by
\[
\left(  T\left(  \sigma\right)  =g_{q}\circ\sigma\circ\left(  g_{p}\right)
^{-1}\ \ \ \ \ \ \ \ \ \ \text{for every }\sigma\in\left\{  \tau\in
S_{n}\ \mid\ \tau\left(  n\right)  =n\right\}  \right)  .
\]
Then, this map $T$ is well-defined and bijective.
\end{lemma}

\begin{proof}
[Proof of Lemma \ref{lem.laplace.gpshort}.]\textbf{(a)} This is trivial.

\textbf{(b)} Let $p\in\left[  n\right]  $. Exercise \ref{exe.perm.cycles}
\textbf{(d)} (applied to $k=p+1$ and $\left(  i_{1},i_{2},\ldots,i_{k}\right)
=\left(  p,p+1,\ldots,n\right)  $) yields%
\[
\left(  -1\right)  ^{\operatorname*{cyc}\nolimits_{p,p+1,\ldots,n}}=\left(
-1\right)  ^{n-\left(  p+1\right)  -1}=\left(  -1\right)  ^{n-p-2}=\left(
-1\right)  ^{n-p}.
\]
Now, $g_{p}=\operatorname*{cyc}\nolimits_{p,p+1,\ldots,n}$, so that $\left(
-1\right)  ^{g_{p}}=\left(  -1\right)  ^{\operatorname*{cyc}%
\nolimits_{p,p+1,\ldots,n}}=\left(  -1\right)  ^{n-p}$. This proves Lemma
\ref{lem.laplace.gpshort} \textbf{(b)}.

\textbf{(c)} We have $g_{p}\left(  n\right)  =p$ (since $g_{p}%
=\operatorname*{cyc}\nolimits_{p,p+1,\ldots,n}$). Also, $g_{p}$ is injective
(since $g_{p}$ is a permutation). Therefore, for every $i\in\left[
n-1\right]  $, we have%
\begin{align*}
g_{p}\left(  i\right)   &  \neq g_{p}\left(  n\right)
\ \ \ \ \ \ \ \ \ \ \left(  \text{since }i\neq n\text{ (because }i\in\left[
n-1\right]  \text{) and since }g_{p}\text{ is injective}\right) \\
&  =p,
\end{align*}
so that $g_{p}\left(  i\right)  \in\left[  n\right]  \setminus\left\{
p\right\}  $. This shows that the map $g_{p}^{\prime}$ is well-defined.

To prove that $g_{p}^{\prime}$ is bijective, we can construct its inverse.
Indeed, for every $i\in\left[  n\right]  \setminus\left\{  p\right\}  $, we
have%
\[
\left(  g_{p}\right)  ^{-1}\left(  i\right)  \neq n\ \ \ \ \ \ \ \ \ \ \left(
\text{since }i\neq p=g_{p}\left(  n\right)  \right)
\]
and thus $\left(  g_{p}\right)  ^{-1}\left(  i\right)  \in\left[  n-1\right]
$. Hence, we can define a map $h:\left[  n\right]  \setminus\left\{
p\right\}  \rightarrow\left[  n-1\right]  $ by%
\[
\left(  h\left(  i\right)  =\left(  g_{p}\right)  ^{-1}\left(  i\right)
\ \ \ \ \ \ \ \ \ \ \text{for every }i\in\left[  n\right]  \setminus\left\{
p\right\}  \right)  .
\]
It is straightforward to check that the maps $g_{p}^{\prime}$ and $h$ are
mutually inverse. Thus, $g_{p}^{\prime}$ is bijective. Lemma
\ref{lem.laplace.gpshort} \textbf{(c)} is thus proven.

\textbf{(d)} We have $g_{p}\left(  n\right)  =p$ (since $g_{p}%
=\operatorname*{cyc}\nolimits_{p,p+1,\ldots,n}$) and $g_{q}\left(  n\right)
=q$ (similarly). Hence, $\left(  g_{p}\right)  ^{-1}\left(  p\right)  =n$
(since $g_{p}\left(  n\right)  =p$) and $\left(  g_{q}\right)  ^{-1}\left(
q\right)  =n$ (since $g_{q}\left(  n\right)  =q$).

For every $\sigma\in\left\{  \tau\in S_{n}\ \mid\ \tau\left(  n\right)
=n\right\}  $, we have $\sigma\left(  n\right)  =n$ and thus
\[
\left(  g_{q}\circ\sigma\circ\left(  g_{p}\right)  ^{-1}\right)  \left(
p\right)  =g_{q}\left(  \sigma\left(  \underbrace{\left(  g_{p}\right)
^{-1}\left(  p\right)  }_{=n}\right)  \right)  =g_{q}\left(
\underbrace{\sigma\left(  n\right)  }_{=n}\right)  =g_{q}\left(  n\right)  =q
\]
and therefore $g_{q}\circ\sigma\circ\left(  g_{p}\right)  ^{-1}\in\left\{
\tau\in S_{n}\ \mid\ \tau\left(  p\right)  =q\right\}  $. Thus, the map $T$ is well-defined.

We can also define a map
\[
Q:\left\{  \tau\in S_{n}\ \mid\ \tau\left(  p\right)  =q\right\}
\rightarrow\left\{  \tau\in S_{n}\ \mid\ \tau\left(  n\right)  =n\right\}
\]
by%
\[
\left(  Q\left(  \sigma\right)  =\left(  g_{q}\right)  ^{-1}\circ\sigma\circ
g_{p}\ \ \ \ \ \ \ \ \ \ \text{for every }\sigma\in\left\{  \tau\in
S_{n}\ \mid\ \tau\left(  p\right)  =q\right\}  \right)  .
\]
The well-definedness of $Q$ can be checked similarly to how we proved the
well-definedness of $T$. It is straightforward to verify that the maps $Q$ and
$T$ are mutually inverse. Thus, $T$ is bijective. This completes the proof of
Lemma \ref{lem.laplace.gpshort} \textbf{(d)}.
\end{proof}
\end{vershort}

\begin{verlong}
\begin{lemma}
\label{lem.laplace.gp}For every $n\in\mathbb{N}$, let $\left[  n\right]  $
denote the set $\left\{  1,2,\ldots,n\right\}  $.

Let $n\in\mathbb{N}$. We recall that, for each $k\in\left\{  1,2,\ldots
,n-1\right\}  $, we have defined $s_{k}$ to be the permutation in $S_{n}$ that
swaps $k$ with $k+1$ but leaves all other numbers unchanged.

For every $p\in\left[  n\right]  $, we define a permutation $g_{p}\in S_{n}$
by%
\[
g_{p}=s_{p}\circ s_{p+1}\circ\cdots\circ s_{n-1}.
\]
(Thus, for $n>0$ and $p=n$, we have $g_{n}=s_{n}\circ s_{n+1}\circ\cdots\circ
s_{n-1}=\left(  \text{a composition of }0\text{ permutations}\right)
=\operatorname*{id}$.)

\textbf{(a)} We have $g_{p}\left(  i\right)  =i$ for every $p\in\left[
n\right]  $ and every $i\in\left[  n\right]  $ satisfying $i<p$.

\textbf{(b)} We have $g_{p}\left(  i\right)  =i+1$ for every $p\in\left[
n\right]  $ and every $i\in\left[  n\right]  $ satisfying $p\leq i<n$.

\textbf{(c)} We have $g_{p}\left(  n\right)  =p$ for every $p\in\left[
n\right]  $.

\textbf{(d)} We have $\left(  g_{p}\left(  1\right)  ,g_{p}\left(  2\right)
,\ldots,g_{p}\left(  n-1\right)  \right)  =\left(  1,2,\ldots,\widehat{p}%
,\ldots,n\right)  $ for every $p\in\left[  n\right]  $.

\textbf{(e)} We have $\left(  -1\right)  ^{g_{p}}=\left(  -1\right)  ^{n-p}$
for every $p\in\left[  n\right]  $.

\textbf{(f)} Let $p\in\left[  n\right]  $. We define a map%
\[
g_{p}^{\prime}:\left[  n-1\right]  \rightarrow\left[  n\right]  \setminus
\left\{  p\right\}
\]
by%
\[
\left(  g_{p}^{\prime}\left(  i\right)  =g_{p}\left(  i\right)
\ \ \ \ \ \ \ \ \ \ \text{for every }i\in\left[  n-1\right]  \right)  .
\]
This map $g_{p}^{\prime}$ is well-defined and bijective.

\textbf{(g)} Let $p\in\left[  n\right]  $ and $q\in\left[  n\right]  $. We
define a map
\[
T:\left\{  \tau\in S_{n}\ \mid\ \tau\left(  n\right)  =n\right\}
\rightarrow\left\{  \tau\in S_{n}\ \mid\ \tau\left(  p\right)  =q\right\}
\]
by
\[
\left(  T\left(  \sigma\right)  =g_{q}\circ\sigma\circ\left(  g_{p}\right)
^{-1}\ \ \ \ \ \ \ \ \ \ \text{for every }\sigma\in\left\{  \tau\in
S_{n}\ \mid\ \tau\left(  n\right)  =n\right\}  \right)  .
\]
Then, this map $T$ is well-defined and bijective.
\end{lemma}

\begin{remark}
Let $n\in\mathbb{N}$. Let us use the notations of Lemma \ref{lem.laplace.gp}
and of Definition \ref{def.perm.cycles}. Then, $g_{p}=\operatorname*{cyc}%
\nolimits_{p,p+1,\ldots,n}$ (as follows from parts \textbf{(a)}, \textbf{(b)}
and \textbf{(c)} of Lemma \ref{lem.laplace.gp}). From this viewpoint, it
appears weird that I have not defined $g_{p}$ as $\operatorname*{cyc}%
\nolimits_{p,p+1,\ldots,n}$. The reason is merely that I wanted to avoid using
cycles (an aesthetical choice).

Parts \textbf{(a)}, \textbf{(b)} and \textbf{(c)} of Lemma
\ref{lem.laplace.gp} can be viewed as an analogue of (\ref{sol.ps2.2.4.c.aik}).
\end{remark}

\begin{proof}
[Proof of Lemma \ref{lem.laplace.gp}.]Let us first recall that if $n>0$, then
$g_{n}$ is well-defined (since $n\in\left[  n\right]  $) and satisfies
\begin{align*}
g_{n}  &  =s_{n}\circ s_{n+1}\circ\cdots\circ s_{n-1}%
\ \ \ \ \ \ \ \ \ \ \left(  \text{by the definition of }g_{n}\right) \\
&  =\left(  \text{a composition of }0\text{ permutations}\right)
=\operatorname*{id}.
\end{align*}
Moreover, every $p\in\left[  n-1\right]  $ satisfies%
\begin{equation}
g_{p}=s_{p}\circ g_{p+1} \label{pf.lem.laplace.gp.recgp}%
\end{equation}
\footnote{\textit{Proof of (\ref{pf.lem.laplace.gp.recgp}):} Let $p\in\left[
n-1\right]  $. Then, $p\in\left[  n-1\right]  =\left\{  1,2,\ldots
,p-1\right\}  $, so that $p+1\in\left\{  2,3,\ldots,p\right\}  \subseteq
\left\{  1,2,\ldots,p\right\}  =\left[  p\right]  $. Hence, $g_{p+1}$ is
well-defined. The definition of $g_{p+1}$ yields $g_{p+1}=s_{p+1}\circ
s_{\left(  p+1\right)  +1}\circ\cdots\circ s_{n-1}=s_{p+1}\circ s_{p+2}%
\circ\cdots\circ s_{n-1}$. Now, $p\in\left[  n-1\right]  \subseteq\left[
n\right]  $, and thus $g_{p}$ is well-defined. The definition of $g_{p}$
yields%
\[
g_{p}=s_{p}\circ s_{p+1}\circ\cdots\circ s_{n-1}=s_{p}\circ\underbrace{\left(
s_{p+1}\circ s_{p+2}\circ\cdots\circ s_{n-1}\right)  }_{=g_{p+1}}=s_{p}\circ
g_{p+1}.
\]
This proves (\ref{pf.lem.laplace.gp.recgp}).}. Moreover, every $p\in\left[
n-1\right]  $ satisfies
\begin{align}
&  s_{p}\left(  p\right)  =p+1;\label{pf.lem.laplace.gp.sp.1}\\
&  s_{p}\left(  p+1\right)  =p;\label{pf.lem.laplace.gp.sp.2}\\
&  \left(  s_{p}\left(  i\right)  =i\ \ \ \ \ \ \ \ \ \ \text{for every }%
i\in\left[  n\right]  \setminus\left\{  p,p+1\right\}  \right)  .
\label{pf.lem.laplace.gp.sp.3}%
\end{align}
\footnote{\textit{Proof of (\ref{pf.lem.laplace.gp.sp.1}),
(\ref{pf.lem.laplace.gp.sp.2}) and (\ref{pf.lem.laplace.gp.sp.3}):} Let
$p\in\left[  n-1\right]  $. Recall that $s_{p}$ is defined as the permutation
in $S_{n}$ that swaps $p$ with $p+1$ but leaves all other numbers unchanged.
Thus, the permutation $s_{p}$ swaps $p$ with $p+1$. In other words, we have
$s_{p}\left(  p\right)  =p+1$ and $s_{p}\left(  p+1\right)  =p$. This proves
(\ref{pf.lem.laplace.gp.sp.1}) and (\ref{pf.lem.laplace.gp.sp.2}).
Furthermore, the permutation $s_{p}$ leaves all other numbers unchanged (where
\textquotedblleft other\textquotedblright\ means \textquotedblleft other than
$p$ and $p+1$\textquotedblright). In other words, $s_{p}\left(  i\right)  =i$
for every $i\in\left[  n\right]  \setminus\left\{  p,p+1\right\}  $. This
proves (\ref{pf.lem.laplace.gp.sp.3}).}

\textbf{(a)} Fix $i\in\left[  n\right]  $. Thus, $i\in\left[  n\right]
=\left\{  1,2,\ldots,n\right\}  $, so that $1\leq i\leq n$. Let us prove that
\begin{equation}
g_{n-q}\left(  i\right)  =i\ \ \ \ \ \ \ \ \ \ \text{for every }q\in\left\{
0,1,\ldots,n-i-1\right\}  . \label{pf.lem.laplace.gp.a.claim}%
\end{equation}

[\textit{Proof of (\ref{pf.lem.laplace.gp.a.claim}):} We shall prove
(\ref{pf.lem.laplace.gp.a.claim}) by induction over $q$:

\textit{Induction base:} If $0\in\left\{  0,1,\ldots,n-i-1\right\}  $, then
$g_{n-0}\left(  i\right)  =i$\ \ \ \ \footnote{\textit{Proof.} Assume that
$0\in\left\{  0,1,\ldots,n-i-1\right\}  $. Thus, $0\leq0\leq n-i-1$, so that
$0\leq n-i-1$ and therefore $1\leq n-\underbrace{i}_{\geq0}\leq n$. Hence,
$n\geq1>0$, and thus $g_{n}=\operatorname*{id}$ (as we know). Hence,
$\underbrace{g_{n-0}}_{=g_{n}=\operatorname*{id}}\left(  i\right)
=\operatorname*{id}\left(  i\right)  =i$, qed.}. In other words,
(\ref{pf.lem.laplace.gp.a.claim}) holds for $q=0$. Thus, the induction base is complete.

\textit{Induction step:} Let $Q\in\left\{  0,1,\ldots,n-i-1\right\}  $ be
positive. Assume that (\ref{pf.lem.laplace.gp.a.claim}) holds for $q=Q-1$. We
need to show that (\ref{pf.lem.laplace.gp.a.claim}) holds for $q=Q$.

We have $Q\leq n-i-1$ (since $Q\in\left\{  0,1,\ldots,n-i-1\right\}  $) and
$Q\geq1$ (since $Q$ is positive and belongs to $\left\{  0,1,\ldots
,n-i-1\right\}  $). Hence, $n-\underbrace{Q}_{\substack{\leq
n-i-1<n-1\\\text{(since }i\geq1>0\text{)}}}>n-\left(  n-1\right)  =1$. Hence,
$n-Q\geq1$. Also, $n-\underbrace{Q}_{\geq1}\leq n-1$. Combining $n-Q\geq1$
with $n-Q\leq n-1$, we obtain $1\leq n-Q\leq n-1$ and thus $n-Q\in\left[
n-1\right]  $. Hence, (\ref{pf.lem.laplace.gp.recgp}) (applied to $p=n-Q$)
shows that $g_{n-Q}=s_{n-Q}\circ g_{n-Q+1}$.

We have $Q\leq n-i-1$, thus $n-\underbrace{Q}_{\leq n-i-1}-1\geq n-\left(
n-i-1\right)  -1=i$, so that $i\leq n-Q-1<n-Q$ and thus $i\neq n-Q$. Also,
$i<n-Q<n-Q+1$ and thus $i\neq n-Q+1$. Combining $i\neq n-Q$ with $i\neq
n-Q+1$, we obtain $i\notin\left\{  n-Q,n-Q+1\right\}  $. Combined with
$i\in\left[  n\right]  $, this yields $i\in\left[  n\right]  \setminus\left\{
n-Q,n-Q+1\right\}  $. Hence, (\ref{pf.lem.laplace.gp.sp.3}) (applied to
$p=n-Q$) shows that $s_{n-Q}\left(  i\right)  =i$.

But we assumed that (\ref{pf.lem.laplace.gp.a.claim}) holds for $q=Q-1$. In
other words, we have $g_{n-\left(  Q-1\right)  }\left(  i\right)  =i$. Since
$n-\left(  Q-1\right)  =n-Q+1$, this rewrites as $g_{n-Q+1}\left(  i\right)
=i$. Now,%
\[
\underbrace{g_{n-Q}}_{=s_{n-Q}\circ g_{n-Q+1}}\left(  i\right)  =\left(
s_{n-Q}\circ g_{n-Q+1}\right)  \left(  i\right)  =s_{n-Q}\left(
\underbrace{g_{n-Q+1}\left(  i\right)  }_{=i}\right)  =s_{n-Q}\left(
i\right)  =i.
\]
In other words, (\ref{pf.lem.laplace.gp.a.claim}) holds for $q=Q$. This
completes the induction step. Thus, the induction proof of
(\ref{pf.lem.laplace.gp.a.claim}) is complete.]

Now, let us forget that we fixed $i$. We thus have proven
(\ref{pf.lem.laplace.gp.a.claim}) for every $i\in\left[  n\right]  $.

Now, let $p\in\left[  n\right]  $ and $i\in\left[  n\right]  $ be such that
$i<p$. Then, $i<p$, so that $p>i$, and thus $p\geq i+1$ (since both $p$ and
$i$ are integers). Hence, $n-\underbrace{p}_{\geq i+1}\leq n-\left(
i+1\right)  =n-i-1$. Also, $p\leq n$ (since $p\in\left[  n\right]  $), so that
$n-p\geq0$. Combining this with $n-p\leq n-i-1$, we obtain $n-p\in\left\{
0,1,\ldots,n-i-1\right\}  $. Hence, we can apply
(\ref{pf.lem.laplace.gp.a.claim}) to $q=n-p$. We thus obtain $g_{n-\left(
n-p\right)  }\left(  i\right)  =i$. Since $n-\left(  n-p\right)  =p$, this
rewrites as $g_{p}\left(  i\right)  =i$. This proves Lemma
\ref{lem.laplace.gp} \textbf{(a)}.

\textbf{(b)} Fix $i\in\left[  n\right]  $ such that $i<n$. Thus, $i\in\left[
n\right]  =\left\{  1,2,\ldots,n\right\}  $, so that $i\geq1$. Combined with
$i<n$, this yields $1\leq i<n$.

We have $1\leq i<n$, so that $i\in\left\{  1,2,\ldots,n-1\right\}  $ and thus
$i+1\in\left\{  2,3,\ldots,n\right\}  \subseteq\left\{  1,2,\ldots,n\right\}
=\left[  n\right]  $.

Let us prove that
\begin{equation}
g_{n-q}\left(  i\right)  =i+1\ \ \ \ \ \ \ \ \ \ \text{for every }q\in\left\{
n-i,n-i+1,\ldots,n-1\right\}  . \label{pf.lem.laplace.gp.b.claim}%
\end{equation}

[\textit{Proof of (\ref{pf.lem.laplace.gp.b.claim}):} We shall prove
(\ref{pf.lem.laplace.gp.b.claim}) by induction over $q$:

\textit{Induction base:} If $n-i\in\left\{  n-i,n-i+1,\ldots,n-1\right\}  $,
then $g_{n-\left(  n-i\right)  }\left(  i\right)  =i+1$%
\ \ \ \ \footnote{\textit{Proof.} Assume that $n-i\in\left\{  n-i,n-i+1,\ldots
,n-1\right\}  $. We have $i+1\in\left[  n\right]  $. Hence, $g_{i+1}$ is
well-defined. Also, $i<i+1$. Hence, Lemma \ref{lem.laplace.gp} \textbf{(a)}
(applied to $p=i+1$) shows that $g_{i+1}\left(  i\right)  =i$.
\par
Now, $i\in\left\{  1,2,\ldots,n-1\right\}  =\left[  n-1\right]  $. Hence,
(\ref{pf.lem.laplace.gp.recgp}) (applied to $p=i$) yields $g_{i}=s_{i}\circ
g_{i+1}$. Hence,%
\[
\underbrace{g_{i}}_{=s_{i}\circ g_{i+1}}\left(  i\right)  =\left(  s_{i}\circ
g_{i+1}\right)  \left(  i\right)  =s_{i}\left(  \underbrace{g_{i+1}\left(
i\right)  }_{=i}\right)  =s_{i}\left(  i\right)  =i+1
\]
(by (\ref{pf.lem.laplace.gp.sp.1}), applied to $p=i$). Since $n-\left(
n-i\right)  =i$, we now have $g_{n-\left(  n-i\right)  }\left(  i\right)
=g_{i}\left(  i\right)  =i+1$. Qed.}. In other words,
(\ref{pf.lem.laplace.gp.b.claim}) holds for $q=n-i$. Thus, the induction base
is complete.

\textit{Induction step:} Let $Q\in\left\{  n-i,n-i+1,\ldots,n-1\right\}  $ be
such that $Q>n-i$. Assume that (\ref{pf.lem.laplace.gp.b.claim}) holds for
$q=Q-1$. We need to show that (\ref{pf.lem.laplace.gp.b.claim}) holds for
$q=Q$.

Now, we have $Q\leq n-1$ (since $Q\in\left\{  n-i,n-i+1,\ldots,n-1\right\}
$). Hence, $n-\underbrace{Q}_{\leq n-1}\geq n-\left(  n-1\right)  =1$. Also,
$n-\underbrace{Q}_{>n-i}<n-\left(  n-i\right)  =i<n$, so that $n-Q\leq n-1$
(since both $n-Q$ and $n$ are integers). Combining $n-Q\geq1$ with $n-Q\leq
n-1$, we obtain $1\leq n-Q\leq n-1$ and thus $n-Q\in\left[  n-1\right]  $.
Hence, (\ref{pf.lem.laplace.gp.recgp}) (applied to $p=n-Q$) shows that
$g_{n-Q}=s_{n-Q}\circ g_{n-Q+1}$.

We have $Q>n-i$, thus $n-\underbrace{Q}_{>n-i}<n-\left(  n-i\right)  =i$.
Hence, $i>n-Q$. Adding $1$ to both sides of this inequality, we obtain
$i+1>n-Q+1$, so that $i+1\neq n-Q+1$. Also, $i+1>n-Q+1>n-Q$, so that $i+1\neq
n-Q$. Combining $i+1\neq n-Q$ with $i+1\neq n-Q+1$, we obtain $i+1\notin%
\left\{  n-Q,n-Q+1\right\}  $. Combined with $i+1\in\left[  n\right]  $, this
yields $i+1\in\left[  n\right]  \setminus\left\{  n-Q,n-Q+1\right\}  $. Hence,
(\ref{pf.lem.laplace.gp.sp.3}) (applied to $n-Q$ and $i+1$ instead of $p$ and
$i$) shows that $s_{n-Q}\left(  i+1\right)  =i+1$.

But we assumed that (\ref{pf.lem.laplace.gp.b.claim}) holds for $q=Q-1$. In
other words, we have $g_{n-\left(  Q-1\right)  }\left(  i\right)  =i+1$. Since
$n-\left(  Q-1\right)  =n-Q+1$, this rewrites as $g_{n-Q+1}\left(  i\right)
=i+1$. Now,%
\[
\underbrace{g_{n-Q}}_{=s_{n-Q}\circ g_{n-Q+1}}\left(  i\right)  =\left(
s_{n-Q}\circ g_{n-Q+1}\right)  \left(  i\right)  =s_{n-Q}\left(
\underbrace{g_{n-Q+1}\left(  i\right)  }_{=i+1}\right)  =s_{n-Q}\left(
i+1\right)  =i+1.
\]
In other words, (\ref{pf.lem.laplace.gp.b.claim}) holds for $q=Q$. This
completes the induction step. Thus, the induction proof of
(\ref{pf.lem.laplace.gp.b.claim}) is complete.]

Now, let us forget that we fixed $i$. We thus have proven
(\ref{pf.lem.laplace.gp.b.claim}) for every $i\in\left[  n\right]  $
satisfying $i<n$.

Now, let $p\in\left[  n\right]  $ and $i\in\left[  n\right]  $ be such that
$p\leq i<n$. Combining $n-\underbrace{p}_{\leq i}\geq n-i$ with
$n-\underbrace{p}_{\substack{\geq1\\\text{(since }p\in\left[  n\right]
\text{)}}}\leq n-1$, we obtain $n-i\leq n-p\leq n-1$, so that $n-p\in\left\{
n-i,n-i+1,\ldots,n-1\right\}  $. Hence, we can apply
(\ref{pf.lem.laplace.gp.b.claim}) to $q=n-p$. We thus obtain $g_{n-\left(
n-p\right)  }\left(  i\right)  =i+1$. Since $n-\left(  n-p\right)  =p$, this
rewrites as $g_{p}\left(  i\right)  =i+1$. This proves Lemma
\ref{lem.laplace.gp} \textbf{(b)}.

\textbf{(c)} Let us first show that%
\begin{equation}
g_{n-q}\left(  n\right)  =n-q\ \ \ \ \ \ \ \ \ \ \text{for every }q\in\left\{
0,1,\ldots,n-1\right\}  . \label{pf.lem.laplace.gp.c.claim}%
\end{equation}

[\textit{Proof of (\ref{pf.lem.laplace.gp.c.claim}):} We shall prove
(\ref{pf.lem.laplace.gp.c.claim}) by induction over $q$:

\textit{Induction base:} If $0\in\left\{  0,1,\ldots,n-1\right\}  $, then
$g_{n-0}\left(  n\right)  =n-0$\ \ \ \ \footnote{\textit{Proof.} Assume that
$0\in\left\{  0,1,\ldots,n-1\right\}  $. Thus, $0\leq0\leq n-1$, so that
$0\leq n-1$ and therefore $1\leq n$. Hence, $n\geq1>0$, and thus
$g_{n}=\operatorname*{id}$ (as we know). Hence, $\underbrace{g_{n-0}}%
_{=g_{n}=\operatorname*{id}}\left(  n\right)  =\operatorname*{id}\left(
n\right)  =n=n-0$, qed.}. In other words, (\ref{pf.lem.laplace.gp.c.claim})
holds for $q=0$. Thus, the induction base is complete.

\textit{Induction step:} Let $Q\in\left\{  0,1,\ldots,n-1\right\}  $ be
positive. Assume that (\ref{pf.lem.laplace.gp.c.claim}) holds for $q=Q-1$. We
need to show that (\ref{pf.lem.laplace.gp.c.claim}) holds for $q=Q$.

We have $Q\leq n-1$ (since $Q\in\left\{  0,1,\ldots,n-1\right\}  $) and
$Q\geq1$ (since $Q$ is positive and belongs to $\left\{  0,1,\ldots
,n-1\right\}  $). Hence, $n-\underbrace{Q}_{\leq n-1}\geq n-\left(
n-1\right)  =1$. Also, $n-\underbrace{Q}_{\geq1}\leq n-1$. Combining
$n-Q\geq1$ with $n-Q\leq n-1$, we obtain $1\leq n-Q\leq n-1$ and thus
$n-Q\in\left[  n-1\right]  $. Hence, (\ref{pf.lem.laplace.gp.recgp}) (applied
to $p=n-Q$) shows that $g_{n-Q}=s_{n-Q}\circ g_{n-Q+1}$.

We have $n-Q\in\left[  n-1\right]  $ (since $1\leq n-Q\leq n-1$). Hence,
$s_{n-Q}\left(  n-Q+1\right)  =n-Q$ (by (\ref{pf.lem.laplace.gp.sp.2}),
applied to $p=n-Q$).

But we assumed that (\ref{pf.lem.laplace.gp.c.claim}) holds for $q=Q-1$. In
other words, we have $g_{n-\left(  Q-1\right)  }\left(  n\right)  =n-\left(
Q-1\right)  $. Since $n-\left(  Q-1\right)  =n-Q+1$, this rewrites as
$g_{n-Q+1}\left(  n\right)  =n-Q+1$. Now,%
\begin{align*}
\underbrace{g_{n-Q}}_{=s_{n-Q}\circ g_{n-Q+1}}\left(  n\right)   &  =\left(
s_{n-Q}\circ g_{n-Q+1}\right)  \left(  n\right)  =s_{n-Q}\left(
\underbrace{g_{n-Q+1}\left(  n\right)  }_{=n-Q+1}\right) \\
&  =s_{n-Q}\left(  n-Q+1\right)  =n-Q.
\end{align*}
In other words, (\ref{pf.lem.laplace.gp.c.claim}) holds for $q=Q$. This
completes the induction step. Thus, the induction proof of
(\ref{pf.lem.laplace.gp.c.claim}) is complete.]

Now, let $p\in\left[  n\right]  $. Then, $p\in\left[  n\right]  =\left\{
1,2,\ldots,n\right\}  $, so that $n-p\in\left\{  0,1,\ldots,n-1\right\}  $.
Thus, we can apply (\ref{pf.lem.laplace.gp.c.claim}) to $q=n-p$. We thus
obtain $g_{n-\left(  n-p\right)  }\left(  n\right)  =n-\left(  n-p\right)  $.
Since $n-\left(  n-p\right)  =p$, this rewrites as $g_{p}\left(  n\right)
=p$. This proves Lemma \ref{lem.laplace.gp} \textbf{(c)}.

\textbf{(d)} In the following argument, we shall use the concept of the
concatenation of two lists. It is defined as follows: If $\mathbf{a}=\left(
a_{1},a_{2},\ldots,a_{k}\right)  $ and $\mathbf{b}=\left(  b_{1},b_{2}%
,\ldots,b_{m}\right)  $ are two finite lists (for example, of integers), then
the \textit{concatenation} of $\mathbf{a}$ with $\mathbf{b}$ is defined to be
the list $\left(  a_{1},a_{2},\ldots,a_{k},b_{1},b_{2},\ldots,b_{m}\right)  $.
(For example, the concatenation of the lists $\left(  1,2,6\right)  $ and
$\left(  0,7\right)  $ is $\left(  1,2,6,0,7\right)  $.)

Let $p\in\left[  n\right]  $. We have $g_{p}\left(  i\right)  =i$ for every
$i\in\left\{  1,2,\ldots,p-1\right\}  $\ \ \ \ \footnote{\textit{Proof.} Let
$i\in\left\{  1,2,\ldots,p-1\right\}  $. Thus, $i\geq1$ and $i\leq p-1$.
Combining $i\geq1$ with $i\leq p-1\leq p\leq n$ (since $p\in\left[  n\right]
$), we obtain $1\leq i\leq n$, so that $i\in\left[  n\right]  $. Also, $i\leq
p-1<p$. Hence, Lemma \ref{lem.laplace.gp} \textbf{(a)} shows that we have
$g_{p}\left(  i\right)  =i$, qed.}. In other words,%
\[
\left(  g_{p}\left(  1\right)  ,g_{p}\left(  2\right)  ,\ldots,g_{p}\left(
p-1\right)  \right)  =\left(  1,2,\ldots,p-1\right)  .
\]

Also, we have $g_{p}\left(  i\right)  =i+1$ for every $i\in\left\{
p,p+1,\ldots,n-1\right\}  $\ \ \ \ \footnote{\textit{Proof.} Let $i\in\left\{
p,p+1,\ldots,n-1\right\}  $. Thus, $i\geq p$ and $i\leq n-1$. Combining $i\geq
p\geq1$ (since $p\in\left[  n\right]  $) with $i\leq n-1\leq n$, we obtain
$1\leq i\leq n$, so that $i\in\left[  n\right]  $. Also, $p\leq i$ (since
$i\geq p$) and $i\leq n-1<n$, so that $p\leq i<n$. Hence, Lemma
\ref{lem.laplace.gp} \textbf{(b)} shows that we have $g_{p}\left(  i\right)
=i+1$, qed.}. In other words,%
\begin{align*}
\left(  g_{p}\left(  p\right)  ,g_{p}\left(  p+1\right)  ,\ldots,g_{p}\left(
n-1\right)  \right)   &  =\left(  p+1,\left(  p+1\right)  +1,\ldots,\left(
n-1\right)  +1\right) \\
&  =\left(  p+1,p+2,\ldots,n\right)  .
\end{align*}

Now,%
\begin{align*}
&  \left(  \text{the concatenation of the list }\underbrace{\left(
g_{p}\left(  1\right)  ,g_{p}\left(  2\right)  ,\ldots,g_{p}\left(
p-1\right)  \right)  }_{=\left(  1,2,\ldots,p-1\right)  }\right. \\
&  \ \ \ \ \ \ \ \ \ \ \left.  \text{with the list }\underbrace{\left(
g_{p}\left(  p\right)  ,g_{p}\left(  p+1\right)  ,\ldots,g_{p}\left(
n-1\right)  \right)  }_{=\left(  p+1,p+2,\ldots,n\right)  }\right) \\
&  =\left(  \text{the concatenation of the list }\left(  1,2,\ldots
,p-1\right)  \text{ with the list }\left(  p+1,p+2,\ldots,n\right)  \right) \\
&  =\left(  1,2,\ldots,p-1,p+1,p+2,\ldots,n\right)  .
\end{align*}
Comparing this with%
\[
\left(  1,2,\ldots,\widehat{p},\ldots,n\right)  =\left(  1,2,\ldots
,p-1,p+1,p+2,\ldots,n\right)  ,
\]
we obtain
\begin{align*}
&  \left(  1,2,\ldots,\widehat{p},\ldots,n\right) \\
&  =\left(  \text{the concatenation of the list }\left(  g_{p}\left(
1\right)  ,g_{p}\left(  2\right)  ,\ldots,g_{p}\left(  p-1\right)  \right)
\right. \\
&  \ \ \ \ \ \ \ \ \ \ \left.  \text{with the list }\left(  g_{p}\left(
p\right)  ,g_{p}\left(  p+1\right)  ,\ldots,g_{p}\left(  n-1\right)  \right)
\right) \\
&  =\left(  g_{p}\left(  1\right)  ,g_{p}\left(  2\right)  ,\ldots
,g_{p}\left(  n-1\right)  \right)  .
\end{align*}
This proves Lemma \ref{lem.laplace.gp} \textbf{(d)}.

\textbf{(e)} Let us first show that%
\begin{equation}
\left(  -1\right)  ^{g_{n-q}}=\left(  -1\right)  ^{q}%
\ \ \ \ \ \ \ \ \ \ \text{for every }q\in\left\{  0,1,\ldots,n-1\right\}  .
\label{pf.lem.laplace.gp.e.claim}%
\end{equation}

[\textit{Proof of (\ref{pf.lem.laplace.gp.e.claim}):} We shall prove
(\ref{pf.lem.laplace.gp.e.claim}) by induction over $q$:

\textit{Induction base:} If $0\in\left\{  0,1,\ldots,n-1\right\}  $, then
$\left(  -1\right)  ^{g_{n-0}}=\left(  -1\right)  ^{0}$%
\ \ \ \ \footnote{\textit{Proof.} Assume that $0\in\left\{  0,1,\ldots
,n-1\right\}  $. Thus, $0\leq0\leq n-1$, so that $0\leq n-1$ and therefore
$1\leq n$. Hence, $n\geq1>0$, and thus $g_{n}=\operatorname*{id}$ (as we
know). Hence, $\left(  -1\right)  ^{g_{n}}=\left(  -1\right)
^{\operatorname*{id}}=1$, so that $\left(  -1\right)  ^{g_{n-0}}=\left(
-1\right)  ^{g_{n}}=1=\left(  -1\right)  ^{0}$, qed.}. In other words,
(\ref{pf.lem.laplace.gp.e.claim}) holds for $q=0$. Thus, the induction base is complete.

\textit{Induction step:} Let $Q\in\left\{  0,1,\ldots,n-1\right\}  $ be
positive. Assume that (\ref{pf.lem.laplace.gp.e.claim}) holds for $q=Q-1$. We
need to show that (\ref{pf.lem.laplace.gp.e.claim}) holds for $q=Q$.

We have $Q\leq n-1$ (since $Q\in\left\{  0,1,\ldots,n-1\right\}  $) and
$Q\geq1$ (since $Q$ is positive and belongs to $\left\{  0,1,\ldots
,n-1\right\}  $). Hence, $n-\underbrace{Q}_{\leq n-1}\geq n-\left(
n-1\right)  =1$. Also, $n-\underbrace{Q}_{\geq1}\leq n-1$. Combining
$n-Q\geq1$ with $n-Q\leq n-1$, we obtain $1\leq n-Q\leq n-1$ and thus
$n-Q\in\left\{  1,2,\ldots,n-1\right\}  =\left[  n-1\right]  $. Hence,
(\ref{pf.lem.laplace.gp.recgp}) (applied to $p=n-Q$) shows that $g_{n-Q}%
=s_{n-Q}\circ g_{n-Q+1}$.

Proposition \ref{prop.perm.signs.basics} \textbf{(b)} shows that $\left(
-1\right)  ^{s_{k}}=-1$ for every $k\in\left\{  1,2,\ldots,n-1\right\}  $.
Applying this to $k=n-Q$, we obtain $\left(  -1\right)  ^{s_{n-Q}}=-1$.

But we assumed that (\ref{pf.lem.laplace.gp.e.claim}) holds for $q=Q-1$. In
other words, we have $\left(  -1\right)  ^{g_{n-\left(  Q-1\right)  }}=\left(
-1\right)  ^{Q-1}$. Since $n-\left(  Q-1\right)  =n-Q+1$, this rewrites as
$\left(  -1\right)  ^{g_{n-Q+1}}=\left(  -1\right)  ^{Q-1}$. Now, from
$g_{n-Q}=s_{n-Q}\circ g_{n-Q+1}$, we obtain%
\begin{align*}
\left(  -1\right)  ^{g_{n-Q}}  &  =\left(  -1\right)  ^{s_{n-Q}\circ
g_{n-Q+1}}=\underbrace{\left(  -1\right)  ^{s_{n-Q}}}_{=-1}\cdot
\underbrace{\left(  -1\right)  ^{g_{n-Q+1}}}_{=\left(  -1\right)  ^{Q-1}}\\
&  \ \ \ \ \ \ \ \ \ \ \left(  \text{by (\ref{eq.sign.prod}), applied to
}\sigma=s_{n-Q}\text{ and }\tau=g_{n-Q+1}\right) \\
&  =\left(  -1\right)  \cdot\left(  -1\right)  ^{Q-1}=\left(  -1\right)
^{\left(  Q-1\right)  +1}=\left(  -1\right)  ^{Q}.
\end{align*}

In other words, (\ref{pf.lem.laplace.gp.e.claim}) holds for $q=Q$. This
completes the induction step. Thus, the induction proof of
(\ref{pf.lem.laplace.gp.e.claim}) is complete.]

Now, let $p\in\left[  n\right]  $. Then, $p\in\left[  n\right]  =\left\{
1,2,\ldots,n\right\}  $, so that $n-p\in\left\{  0,1,\ldots,n-1\right\}  $.
Thus, we can apply (\ref{pf.lem.laplace.gp.e.claim}) to $q=n-p$. We thus
obtain $\left(  -1\right)  ^{g_{n-\left(  n-p\right)  }}=\left(  -1\right)
^{n-p}$. Since $n-\left(  n-p\right)  =p$, this rewrites as $\left(
-1\right)  ^{g_{p}}=\left(  -1\right)  ^{n-p}$. This proves Lemma
\ref{lem.laplace.gp} \textbf{(e)}.

\textbf{(f)} We have $g_{p}\in S_{n}$. In other words, $g_{p}$ is a
permutation of $\left\{  1,2,\ldots,n\right\}  $ (since $S_{n}$ is the set of
all permutations of $\left\{  1,2,\ldots,n\right\}  $). In other words,
$g_{p}$ is a bijective map $\left\{  1,2,\ldots,n\right\}  \rightarrow\left\{
1,2,\ldots,n\right\}  $. Thus, the map $g_{p}$ is injective and surjective.

For every $i\in\left[  n-1\right]  $, we have $g_{p}\left(  i\right)
\in\left[  n\right]  \setminus\left\{  p\right\}  $%
\ \ \ \ \footnote{\textit{Proof.} Let $i\in\left[  n-1\right]  $. Then,
$i\in\left[  n-1\right]  =\left\{  1,2,\ldots,n-1\right\}  $, so that $i\leq
n-1<n$ and thus $i\neq n$.
\par
Also, $i\in\left\{  1,2,\ldots,n-1\right\}  \subseteq\left\{  1,2,\ldots
,n\right\}  $. Hence, $g_{p}\left(  i\right)  $ is well-defined. Moreover,
$i\in\left\{  1,2,\ldots,n-1\right\}  $, so that $i\geq1$.
\par
The map $g_{p}$ is injective. Since $i\neq n$, we therefore have $g_{p}\left(
i\right)  \neq g_{p}\left(  n\right)  =p$ (by Lemma \ref{lem.laplace.gp}
\textbf{(c)}).
\par
Combining $g_{p}\left(  i\right)  \in\left\{  1,2,\ldots,n\right\}  =\left[
n\right]  $ with $g_{p}\left(  i\right)  \neq p$, we obtain $g_{p}\left(
i\right)  \in\left[  n\right]  \setminus\left\{  p\right\}  $, qed.}. Thus,
the map $g_{p}^{\prime}$ is well-defined. It remains to show that this map
$g_{p}^{\prime}$ is bijective.

For every $j\in\left[  n\right]  \setminus\left\{  p\right\}  $, we have
$\left(  g_{p}\right)  ^{-1}\left(  j\right)  \in\left[  n-1\right]
$\ \ \ \ \footnote{\textit{Proof.} Let $j\in\left[  n\right]  \setminus
\left\{  p\right\}  $. Thus, $j\in\left[  n\right]  $ and $j\neq p$.
\par
Let $i=\left(  g_{p}\right)  ^{-1}\left(  j\right)  $. (This is clearly
well-defined, since $j\in\left[  n\right]  $.) Now, $g_{p}\left(  i\right)
=j$ (since $i=\left(  g_{p}\right)  ^{-1}\left(  j\right)  $). If we had
$i=n$, then we would have $g_{p}\left(  \underbrace{i}_{=n}\right)
=g_{p}\left(  n\right)  =p$ (by Lemma \ref{lem.laplace.gp} \textbf{(c)}),
which would contradict $g_{p}\left(  i\right)  =j\neq p$. Hence, we cannot
have $i=n$. We thus have $i\neq n$.
\par
But $i=\left(  g_{p}\right)  ^{-1}\left(  j\right)  \in\left[  n\right]
=\left\{  1,2,\ldots,n\right\}  $. Combining this with $i\neq n$, we obtain
$i\in\left\{  1,2,\ldots,n\right\}  \setminus\left\{  n\right\}  =\left\{
1,2,\ldots,n-1\right\}  =\left[  n-1\right]  $. Thus, $\left(  g_{p}\right)
^{-1}\left(  j\right)  =i\in\left[  n-1\right]  $, qed.}. Therefore, we can
define a map $h:\left[  n\right]  \setminus\left\{  p\right\}  \rightarrow
\left[  n-1\right]  $ by%
\[
\left(  h\left(  j\right)  =\left(  g_{p}\right)  ^{-1}\left(  j\right)
\ \ \ \ \ \ \ \ \ \ \text{for every }j\in\left[  n\right]  \setminus\left\{
p\right\}  \right)  .
\]
Consider this map $h$.

We have $h\circ g_{p}^{\prime}=\operatorname*{id}$%
\ \ \ \ \footnote{\textit{Proof.} For every $i\in\left[  n-1\right]  $, we
have%
\begin{align*}
\left(  h\circ g_{p}^{\prime}\right)  \left(  i\right)   &  =h\left(
\underbrace{g_{p}^{\prime}\left(  i\right)  }_{\substack{=g_{p}\left(
i\right)  \\\text{(by the definition of }g_{p}^{\prime}\text{)}}}\right)
=h\left(  g_{p}\left(  i\right)  \right)  =\left(  g_{p}\right)  ^{-1}\left(
g_{p}\left(  i\right)  \right)  \ \ \ \ \ \ \ \ \ \ \left(  \text{by the
definition of }h\right) \\
&  =i=\operatorname*{id}\left(  i\right)  .
\end{align*}
Thus, $h\circ g_{p}^{\prime}=\operatorname*{id}$, qed.} and $g_{p}^{\prime
}\circ h=\operatorname*{id}$\ \ \ \ \footnote{\textit{Proof.} For every
$j\in\left[  n\right]  \setminus\left\{  p\right\}  $, we have%
\begin{align*}
\left(  g_{p}^{\prime}\circ h\right)  \left(  j\right)   &  =g_{p}^{\prime
}\left(  \underbrace{h\left(  j\right)  }_{\substack{=\left(  g_{p}\right)
^{-1}\left(  j\right)  \\\text{(by the definition of }h\text{)}}}\right) \\
&  =g_{p}^{\prime}\left(  \left(  g_{p}\right)  ^{-1}\left(  j\right)
\right)  =g_{p}\left(  \left(  g_{p}\right)  ^{-1}\left(  j\right)  \right)
\ \ \ \ \ \ \ \ \ \ \left(  \text{by the definition of }g_{p}^{\prime}\right)
\\
&  =j=\operatorname*{id}\left(  j\right)  .
\end{align*}
Thus, $g_{p}^{\prime}\circ h=\operatorname*{id}$, qed.}. Hence, the maps $h$
and $g_{p}^{\prime}$ are mutually inverse. Thus, the map $g_{p}^{\prime}$ is
invertible. In other words, the map $g_{p}^{\prime}$ is bijective. This
completes the proof of Lemma \ref{lem.laplace.gp} \textbf{(f)}.

\textbf{(g)} For every $\sigma\in\left\{  \tau\in S_{n}\ \mid\ \tau\left(
n\right)  =n\right\}  $, we have%
\[
g_{q}\circ\sigma\circ\left(  g_{p}\right)  ^{-1}\in\left\{  \tau\in
S_{n}\ \mid\ \tau\left(  p\right)  =q\right\}
\]
\footnote{\textit{Proof.} Let $\sigma\in\left\{  \tau\in S_{n}\ \mid
\ \tau\left(  n\right)  =n\right\}  $. Thus, $\sigma$ is an element $\tau$ of
$S_{n}$ satisfying $\tau\left(  n\right)  =n$. In other words, $\sigma$ is an
element of $S_{n}$ and satisfies $\sigma\left(  n\right)  =n$.
\par
We have $g_{p}\left(  n\right)  =p$ (by Lemma \ref{lem.laplace.gp}
\textbf{(c)}). Thus, $\left(  g_{p}\right)  ^{-1}\left(  p\right)  =n$.
Moreover, $g_{q}\left(  n\right)  =q$ (by Lemma \ref{lem.laplace.gp}
\textbf{(c)}, applied to $q$ instead of $p$). Now,%
\[
\left(  g_{q}\circ\sigma\circ\left(  g_{p}\right)  ^{-1}\right)  \left(
p\right)  =g_{q}\left(  \sigma\left(  \underbrace{\left(  g_{p}\right)
^{-1}\left(  p\right)  }_{=n}\right)  \right)  =g_{q}\left(
\underbrace{\sigma\left(  n\right)  }_{=n}\right)  =g_{q}\left(  n\right)
=q.
\]
\par
Now, we know that $g_{q}\circ\sigma\circ\left(  g_{p}\right)  ^{-1}\in S_{n}$
(since $g_{q}$, $\sigma$ and $\left(  g_{p}\right)  ^{-1}$ all belong to
$S_{n}$) and satisfies $\left(  g_{q}\circ\sigma\circ\left(  g_{p}\right)
^{-1}\right)  \left(  p\right)  =q$. In other words, $g_{q}\circ\sigma
\circ\left(  g_{p}\right)  ^{-1}$ is an element $\tau$ of $S_{n}$ satisfying
$\tau\left(  p\right)  =q$. In other words,%
\[
g_{q}\circ\sigma\circ\left(  g_{p}\right)  ^{-1}\in\left\{  \tau\in
S_{n}\ \mid\ \tau\left(  p\right)  =q\right\}  ,
\]
qed.}. Thus, the map $T$ is well-defined.

Furthermore, for every $\sigma\in\left\{  \tau\in S_{n}\ \mid\ \tau\left(
p\right)  =q\right\}  $, we have
\[
\left(  g_{q}\right)  ^{-1}\circ\sigma\circ g_{p}\in\left\{  \tau\in
S_{n}\ \mid\ \tau\left(  n\right)  =n\right\}
\]
\footnote{\textit{Proof.} Let $\sigma\in\left\{  \tau\in S_{n}\ \mid
\ \tau\left(  p\right)  =q\right\}  $. Thus, $\sigma$ is an element $\tau$ of
$S_{n}$ satisfying $\tau\left(  p\right)  =q$. In other words, $\sigma$ is an
element of $S_{n}$ and satisfies $\sigma\left(  p\right)  =q$.
\par
We have $g_{p}\left(  n\right)  =p$ (by Lemma \ref{lem.laplace.gp}
\textbf{(c)}). Moreover, $g_{q}\left(  n\right)  =q$ (by Lemma
\ref{lem.laplace.gp} \textbf{(c)}, applied to $q$ instead of $p$), and thus
$\left(  g_{q}\right)  ^{-1}\left(  q\right)  =n$. Now,%
\[
\left(  \left(  g_{q}\right)  ^{-1}\circ\sigma\circ g_{p}\right)  \left(
n\right)  =\left(  g_{q}\right)  ^{-1}\left(  \sigma\left(  \underbrace{g_{p}%
\left(  n\right)  }_{=p}\right)  \right)  =\left(  g_{q}\right)  ^{-1}\left(
\underbrace{\sigma\left(  p\right)  }_{=q}\right)  =\left(  g_{q}\right)
^{-1}\left(  q\right)  =n.
\]
\par
Now, we know that $\left(  g_{q}\right)  ^{-1}\circ\sigma\circ g_{p}\in S_{n}$
(since $\left(  g_{q}\right)  ^{-1}$, $\sigma$ and $g_{p}$ all belong to
$S_{n}$) and satisfies $\left(  \left(  g_{q}\right)  ^{-1}\circ\sigma\circ
g_{p}\right)  \left(  n\right)  =n$. In other words, $\left(  g_{q}\right)
^{-1}\circ\sigma\circ g_{p}$ is an element $\tau$ of $S_{n}$ satisfying
$\tau\left(  p\right)  =n$. In other words,%
\[
\left(  g_{q}\right)  ^{-1}\circ\sigma\circ g_{p}\in\left\{  \tau\in
S_{n}\ \mid\ \tau\left(  n\right)  =n\right\}  ,
\]
qed.}. Hence, we can define a map
\[
Q:\left\{  \tau\in S_{n}\ \mid\ \tau\left(  p\right)  =q\right\}
\rightarrow\left\{  \tau\in S_{n}\ \mid\ \tau\left(  n\right)  =n\right\}
\]
by%
\[
\left(  Q\left(  \sigma\right)  =\left(  g_{q}\right)  ^{-1}\circ\sigma\circ
g_{p}\ \ \ \ \ \ \ \ \ \ \text{for every }\sigma\in\left\{  \tau\in
S_{n}\ \mid\ \tau\left(  p\right)  =q\right\}  \right)  .
\]
Consider this map $Q$.

We have $T\circ Q=\operatorname*{id}$\ \ \ \ \footnote{\textit{Proof.} Every
$\sigma\in\left\{  \tau\in S_{n}\ \mid\ \tau\left(  p\right)  =q\right\}  $
satisfies%
\begin{align*}
\left(  T\circ Q\right)  \left(  \sigma\right)   &  =T\left(
\underbrace{Q\left(  \sigma\right)  }_{\substack{=\left(  g_{q}\right)
^{-1}\circ\sigma\circ g_{p}\\\text{(by the definition of }Q\text{)}}}\right)
=T\left(  \left(  g_{q}\right)  ^{-1}\circ\sigma\circ g_{p}\right) \\
&  =\underbrace{g_{q}\circ\left(  g_{q}\right)  ^{-1}}_{=\operatorname*{id}%
}\circ\sigma\circ\underbrace{g_{p}\circ\left(  g_{p}\right)  ^{-1}%
}_{=\operatorname*{id}}\ \ \ \ \ \ \ \ \ \ \left(  \text{by the definition of
}T\right) \\
&  =\sigma=\operatorname*{id}\left(  \sigma\right)  .
\end{align*}
In other words, $T\circ Q=\operatorname*{id}$, qed.} and $Q\circ
T=\operatorname*{id}$\ \ \ \ \footnote{\textit{Proof.} Every $\sigma
\in\left\{  \tau\in S_{n}\ \mid\ \tau\left(  n\right)  =n\right\}  $ satisfies%
\begin{align*}
\left(  Q\circ T\right)  \left(  \sigma\right)   &  =Q\left(
\underbrace{T\left(  \sigma\right)  }_{\substack{=g_{q}\circ\sigma\circ\left(
g_{p}\right)  ^{-1}\\\text{(by the definition of }T\text{)}}}\right)
=Q\left(  g_{q}\circ\sigma\circ\left(  g_{p}\right)  ^{-1}\right) \\
&  =\underbrace{\left(  g_{q}\right)  ^{-1}\circ\left(  g_{q}\right)
}_{=\operatorname*{id}}\circ\sigma\circ\underbrace{\left(  g_{p}\right)
^{-1}\circ g_{p}}_{=\operatorname*{id}}\ \ \ \ \ \ \ \ \ \ \left(  \text{by
the definition of }Q\right) \\
&  =\sigma=\operatorname*{id}\left(  \sigma\right)  .
\end{align*}
In other words, $Q\circ T=\operatorname*{id}$, qed.}. Hence, the maps $T$ and
$Q$ are mutually inverse. Thus, the map $T$ is invertible. In other words, the
map $T$ is bijective. This completes the proof of Lemma \ref{lem.laplace.gp}
\textbf{(g)}.
\end{proof}
\end{verlong}

Our next step towards the proof of Theorem \ref{thm.laplace.gen} is the
following lemma:

\begin{lemma}
\label{lem.laplace.Apq}Let $n\in\mathbb{N}$. Let $A=\left(  a_{i,j}\right)
_{1\leq i\leq n,\ 1\leq j\leq n}$ be an $n\times n$-matrix. Let $p\in\left\{
1,2,\ldots,n\right\}  $ and $q\in\left\{  1,2,\ldots,n\right\}  $. Then,%
\[
\sum_{\substack{\sigma\in S_{n};\\\sigma\left(  p\right)  =q}}\left(
-1\right)  ^{\sigma}\prod_{\substack{i\in\left\{  1,2,\ldots,n\right\}
;\\i\neq p}}a_{i,\sigma\left(  i\right)  }=\left(  -1\right)  ^{p+q}%
\det\left(  A_{\sim p,\sim q}\right)  .
\]

\end{lemma}

\begin{vershort}
\begin{proof}
[Proof of Lemma \ref{lem.laplace.Apq}.]Let us use all notations introduced in
Lemma \ref{lem.laplace.gpshort}.

We have $p\in\left\{  1,2,\ldots,n\right\}  =\left[  n\right]  $. Hence,
$g_{p}$ is well-defined. Similarly, $g_{q}$ is well-defined. We have%
\begin{equation}
\left(  g_{p}\left(  1\right)  ,g_{p}\left(  2\right)  ,\ldots,g_{p}\left(
n-1\right)  \right)  =\left(  1,2,\ldots,\widehat{p},\ldots,n\right)
\label{pf.lem.laplace.Apq.short.indices1}%
\end{equation}
(by Lemma \ref{lem.laplace.gpshort} \textbf{(a)}) and
\begin{equation}
\left(  g_{q}\left(  1\right)  ,g_{q}\left(  2\right)  ,\ldots,g_{q}\left(
n-1\right)  \right)  =\left(  1,2,\ldots,\widehat{q},\ldots,n\right)
\label{pf.lem.laplace.Apq.short.indices2}%
\end{equation}
(by Lemma \ref{lem.laplace.gpshort} \textbf{(a)}, applied to $q$ instead of
$p$). Now, the definition of $A_{\sim p,\sim q}$ yields%
\begin{align}
A_{\sim p,\sim q}  &  =\operatorname*{sub}\nolimits_{1,2,\ldots,\widehat{p}%
,\ldots,n}^{1,2,\ldots,\widehat{q},\ldots,n}A=\operatorname*{sub}%
\nolimits_{g_{p}\left(  1\right)  ,g_{p}\left(  2\right)  ,\ldots,g_{p}\left(
n-1\right)  }^{g_{q}\left(  1\right)  ,g_{q}\left(  2\right)  ,\ldots
,g_{q}\left(  n-1\right)  }A\ \ \ \ \ \ \ \ \ \ \left(  \text{by
(\ref{pf.lem.laplace.Apq.short.indices1}) and
(\ref{pf.lem.laplace.Apq.short.indices2})}\right) \nonumber\\
&  =\left(  a_{g_{p}\left(  x\right)  ,g_{q}\left(  y\right)  }\right)
_{1\leq x\leq n-1,\ 1\leq y\leq n-1}\nonumber\\
&  \ \ \ \ \ \ \ \ \ \ \left(  \text{by the definition of }\operatorname*{sub}%
\nolimits_{g_{p}\left(  1\right)  ,g_{p}\left(  2\right)  ,\ldots,g_{p}\left(
n-1\right)  }^{g_{q}\left(  1\right)  ,g_{q}\left(  2\right)  ,\ldots
,g_{q}\left(  n-1\right)  }A\right) \nonumber\\
&  =\left(  a_{g_{p}\left(  i\right)  ,g_{q}\left(  j\right)  }\right)
_{1\leq i\leq n-1,\ 1\leq j\leq n-1}\label{pf.lem.laplace.Apq.short.A}\\
&  \ \ \ \ \ \ \ \ \ \ \left(  \text{here, we renamed the index }\left(
x,y\right)  \text{ as }\left(  i,j\right)  \right)  .\nonumber
\end{align}

Also, $\left[  n\right]  $ is nonempty (since $p\in\left[  n\right]  $), and
thus we have $n>0$.

Now, let us recall the map $T:\left\{  \tau\in S_{n}\ \mid\ \tau\left(
n\right)  =n\right\}  \rightarrow\left\{  \tau\in S_{n}\ \mid\ \tau\left(
p\right)  =q\right\}  $ defined in Lemma \ref{lem.laplace.gpshort}
\textbf{(d)}. Lemma \ref{lem.laplace.gpshort} \textbf{(d)} says that this map
$T$ is well-defined and bijective. Every $\sigma\in\left\{  \tau\in
S_{n}\ \mid\ \tau\left(  n\right)  =n\right\}  $ satisfies%
\begin{equation}
\left(  -1\right)  ^{T\left(  \sigma\right)  }=\left(  -1\right)  ^{p+q}%
\cdot\left(  -1\right)  ^{\sigma} \label{pf.lem.laplace.Apq.short.1}%
\end{equation}
\footnote{\textit{Proof of (\ref{pf.lem.laplace.Apq.short.1}):} Let $\sigma
\in\left\{  \tau\in S_{n}\ \mid\ \tau\left(  n\right)  =n\right\}  $. Applying
Lemma \ref{lem.laplace.gpshort} \textbf{(b)} to $q$ instead of $p$, we obtain
$\left(  -1\right)  ^{g_{q}}=\left(  -1\right)  ^{n-q}=\left(  -1\right)
^{n+q}$ (since $n-q\equiv n+q\operatorname{mod}2$).
\par
The definition of $T\left(  \sigma\right)  $ yields $T\left(  \sigma\right)
=g_{q}\circ\sigma\circ\left(  g_{p}\right)  ^{-1}$. Thus,%
\[
\underbrace{T\left(  \sigma\right)  }_{=g_{q}\circ\sigma\circ\left(
g_{p}\right)  ^{-1}}\circ g_{p}=g_{q}\circ\sigma\circ\underbrace{\left(
g_{p}\right)  ^{-1}\circ g_{p}}_{=\operatorname*{id}}=g_{q}\circ\sigma,
\]
so that%
\begin{align*}
\left(  -1\right)  ^{T\left(  \sigma\right)  \circ g_{p}}  &  =\left(
-1\right)  ^{g_{q}\circ\sigma}=\underbrace{\left(  -1\right)  ^{g_{q}}%
}_{=\left(  -1\right)  ^{n+q}}\cdot\left(  -1\right)  ^{\sigma}%
\ \ \ \ \ \ \ \ \ \ \left(  \text{by (\ref{eq.sign.prod}), applied to }%
g_{q}\text{ and }\sigma\text{ instead of }\sigma\text{ and }\tau\right) \\
&  =\left(  -1\right)  ^{n+q}\cdot\left(  -1\right)  ^{\sigma}.
\end{align*}
Compared with%
\begin{align*}
\left(  -1\right)  ^{T\left(  \sigma\right)  \circ g_{p}}  &  =\left(
-1\right)  ^{T\left(  \sigma\right)  }\cdot\underbrace{\left(  -1\right)
^{g_{p}}}_{\substack{=\left(  -1\right)  ^{n-p}\\\text{(by Lemma
\ref{lem.laplace.gpshort} \textbf{(b)})}}}\ \ \ \ \ \ \ \ \ \ \left(  \text{by
(\ref{eq.sign.prod}), applied to }T\left(  \sigma\right)  \text{ and }%
g_{p}\text{ instead of }\sigma\text{ and }\tau\right) \\
&  =\left(  -1\right)  ^{T\left(  \sigma\right)  }\cdot\left(  -1\right)
^{n-p},
\end{align*}
this yields
\[
\left(  -1\right)  ^{T\left(  \sigma\right)  }\cdot\left(  -1\right)
^{n-p}=\left(  -1\right)  ^{n+q}\cdot\left(  -1\right)  ^{\sigma}.
\]
We can divide both sides of this equality by $\left(  -1\right)  ^{n-p}$
(since $\left(  -1\right)  ^{n-p}\in\left\{  1,-1\right\}  $ is clearly an
invertible integer), and thus we obtain%
\[
\left(  -1\right)  ^{T\left(  \sigma\right)  }=\dfrac{\left(  -1\right)
^{n+q}\cdot\left(  -1\right)  ^{\sigma}}{\left(  -1\right)  ^{n-p}%
}=\underbrace{\dfrac{\left(  -1\right)  ^{n+q}}{\left(  -1\right)  ^{n-p}}%
}_{\substack{=\left(  -1\right)  ^{\left(  n+q\right)  -\left(  n-p\right)
}=\left(  -1\right)  ^{p+q}\\\text{(since }\left(  n+q\right)  -\left(
n-p\right)  =p+q\text{)}}}\cdot\left(  -1\right)  ^{\sigma}=\left(  -1\right)
^{p+q}\cdot\left(  -1\right)  ^{\sigma}.
\]
This proves (\ref{pf.lem.laplace.Apq.short.1}).} and%
\begin{equation}
\prod_{\substack{i\in\left\{  1,2,\ldots,n\right\}  ;\\i\neq p}}a_{i,\left(
T\left(  \sigma\right)  \right)  \left(  i\right)  }=\prod_{i=1}^{n-1}%
a_{g_{p}\left(  i\right)  ,g_{q}\left(  \sigma\left(  i\right)  \right)  }
\label{pf.lem.laplace.Apq.short.2}%
\end{equation}
\footnote{\textit{Proof of (\ref{pf.lem.laplace.Apq.short.2}):} Let $\sigma
\in\left\{  \tau\in S_{n}\ \mid\ \tau\left(  n\right)  =n\right\}  $. Let us
recall the map $g_{p}^{\prime}:\left[  n-1\right]  \rightarrow\left[
n\right]  \setminus\left\{  p\right\}  $ introduced in Lemma
\ref{lem.laplace.gpshort} \textbf{(c)}. Lemma \ref{lem.laplace.gpshort}
\textbf{(c)} says that this map $g_{p}^{\prime}$ is well-defined and
bijective. In other words, $g_{p}^{\prime}$ is a bijection.
\par
Let $i\in\left[  n-1\right]  $. Then, $g_{p}^{\prime}\left(  i\right)
=g_{p}\left(  i\right)  $ (by the definition of $g_{p}^{\prime}$). Also, the
definition of $T$ yields $T\left(  \sigma\right)  =g_{q}\circ\sigma
\circ\left(  g_{p}\right)  ^{-1}$, so that%
\[
\left(  \underbrace{T\left(  \sigma\right)  }_{=g_{q}\circ\sigma\circ\left(
g_{p}\right)  ^{-1}}\right)  \left(  \underbrace{g_{p}^{\prime}\left(
i\right)  }_{=g_{p}\left(  i\right)  }\right)  =\left(  g_{q}\circ\sigma
\circ\left(  g_{p}\right)  ^{-1}\right)  \left(  g_{p}\left(  i\right)
\right)  =g_{q}\left(  \sigma\left(  \underbrace{\left(  g_{p}\right)
^{-1}\left(  g_{p}\left(  i\right)  \right)  }_{=i}\right)  \right)
=g_{q}\left(  \sigma\left(  i\right)  \right)  .
\]
\par
From $g_{p}^{\prime}\left(  i\right)  =g_{p}\left(  i\right)  $ and $\left(
T\left(  \sigma\right)  \right)  \left(  g_{p}^{\prime}\left(  i\right)
\right)  =g_{q}\left(  \sigma\left(  i\right)  \right)  $, we obtain%
\begin{equation}
a_{g_{p}^{\prime}\left(  i\right)  ,\left(  T\left(  \sigma\right)  \right)
\left(  g_{p}^{\prime}\left(  i\right)  \right)  }=a_{g_{p}\left(  i\right)
,g_{q}\left(  \sigma\left(  i\right)  \right)  }.
\label{pf.lem.laplace.Apq.short.2.pf.1}%
\end{equation}
\par
Now, let us forget that we fixed $i$. We thus have proven
(\ref{pf.lem.laplace.Apq.short.2.pf.1}) for every $i\in\left[  n-1\right]  $.
But now, we have%
\begin{align*}
&  \underbrace{\prod_{\substack{i\in\left\{  1,2,\ldots,n\right\}  ;\\i\neq
p}}}_{\substack{=\prod_{\substack{i\in\left[  n\right]  ;\\i\neq
p}}\\\text{(since }\left\{  1,2,\ldots,n\right\}  =\left[  n\right]  \text{)}%
}}a_{i,\left(  T\left(  \sigma\right)  \right)  \left(  i\right)  }\\
&  =\underbrace{\prod_{\substack{i\in\left[  n\right]  ;\\i\neq p}}}%
_{=\prod_{i\in\left[  n\right]  \setminus\left\{  p\right\}  }}a_{i,\left(
T\left(  \sigma\right)  \right)  \left(  i\right)  }=\prod_{i\in\left[
n\right]  \setminus\left\{  p\right\}  }a_{i,\left(  T\left(  \sigma\right)
\right)  \left(  i\right)  }=\underbrace{\prod_{i\in\left[  n-1\right]  }%
}_{=\prod_{i=1}^{n-1}}\underbrace{a_{g_{p}^{\prime}\left(  i\right)  ,\left(
T\left(  \sigma\right)  \right)  \left(  g_{p}^{\prime}\left(  i\right)
\right)  }}_{\substack{=a_{g_{p}\left(  i\right)  ,g_{q}\left(  \sigma\left(
i\right)  \right)  }\\\text{(by (\ref{pf.lem.laplace.Apq.short.2.pf.1}))}}}\\
&  \ \ \ \ \ \ \ \ \ \ \left(
\begin{array}
[c]{c}%
\text{here, we have substituted }g_{p}^{\prime}\left(  i\right)  \text{ for
}i\text{, since}\\
g_{p}^{\prime}:\left[  n-1\right]  \rightarrow\left[  n\right]  \setminus
\left\{  p\right\}  \text{ is a bijection}%
\end{array}
\right) \\
&  =\prod_{i=1}^{n-1}a_{g_{p}\left(  i\right)  ,g_{q}\left(  \sigma\left(
i\right)  \right)  }.
\end{align*}
This proves (\ref{pf.lem.laplace.Apq.short.2}).}.

Now,
\begin{align*}
&  \underbrace{\sum_{\substack{\sigma\in S_{n};\\\sigma\left(  p\right)  =q}%
}}_{=\sum_{\sigma\in\left\{  \tau\in S_{n}\ \mid\ \tau\left(  p\right)
=q\right\}  }}\left(  -1\right)  ^{\sigma}\prod_{\substack{i\in\left\{
1,2,\ldots,n\right\}  ;\\i\neq p}}a_{i,\sigma\left(  i\right)  }\\
&  =\sum_{\sigma\in\left\{  \tau\in S_{n}\ \mid\ \tau\left(  p\right)
=q\right\}  }\left(  -1\right)  ^{\sigma}\prod_{\substack{i\in\left\{
1,2,\ldots,n\right\}  ;\\i\neq p}}a_{i,\sigma\left(  i\right)  }\\
&  =\underbrace{\sum_{\sigma\in\left\{  \tau\in S_{n}\ \mid\ \tau\left(
n\right)  =n\right\}  }}_{=\sum_{\substack{\sigma\in S_{n};\\\sigma\left(
n\right)  =n}}}\underbrace{\left(  -1\right)  ^{T\left(  \sigma\right)  }%
}_{\substack{=\left(  -1\right)  ^{p+q}\cdot\left(  -1\right)  ^{\sigma
}\\\text{(by (\ref{pf.lem.laplace.Apq.short.1}))}}}\underbrace{\prod
_{\substack{i\in\left\{  1,2,\ldots,n\right\}  ;\\i\neq p}}a_{i,\left(
T\left(  \sigma\right)  \right)  \left(  i\right)  }}_{\substack{=\prod
_{i=1}^{n-1}a_{g_{p}\left(  i\right)  ,g_{q}\left(  \sigma\left(  i\right)
\right)  }\\\text{(by (\ref{pf.lem.laplace.Apq.short.2}))}}}\\
&  \ \ \ \ \ \ \ \ \ \ \left(
\begin{array}
[c]{c}%
\text{here, we have substituted }T\left(  \sigma\right)  \text{ for }%
\sigma\text{ in the sum,}\\
\text{since the map }T:\left\{  \tau\in S_{n}\ \mid\ \tau\left(  n\right)
=n\right\}  \rightarrow\left\{  \tau\in S_{n}\ \mid\ \tau\left(  p\right)
=q\right\} \\
\text{is a bijection}%
\end{array}
\right) \\
&  =\sum_{\substack{\sigma\in S_{n};\\\sigma\left(  n\right)  =n}}\left(
-1\right)  ^{p+q}\cdot\left(  -1\right)  ^{\sigma}\prod_{i=1}^{n-1}%
a_{g_{p}\left(  i\right)  ,g_{q}\left(  \sigma\left(  i\right)  \right)  }\\
&  =\left(  -1\right)  ^{p+q}\underbrace{\sum_{\substack{\sigma\in
S_{n};\\\sigma\left(  n\right)  =n}}\left(  -1\right)  ^{\sigma}\prod
_{i=1}^{n-1}a_{g_{p}\left(  i\right)  ,g_{q}\left(  \sigma\left(  i\right)
\right)  }}_{\substack{=\det\left(  \left(  a_{g_{p}\left(  i\right)
,g_{q}\left(  j\right)  }\right)  _{1\leq i\leq n-1,\ 1\leq j\leq n-1}\right)
\\\text{(by Lemma \ref{lem.laplace.lem}, applied to }a_{g_{p}\left(  i\right)
,g_{q}\left(  j\right)  }\text{ instead of }a_{i,j}\text{)}}}\\
&  =\left(  -1\right)  ^{p+q}\det\left(  \underbrace{\left(  a_{g_{p}\left(
i\right)  ,g_{q}\left(  j\right)  }\right)  _{1\leq i\leq n-1,\ 1\leq j\leq
n-1}}_{\substack{=A_{\sim p,\sim q}\\\text{(by
(\ref{pf.lem.laplace.Apq.short.A}))}}}\right)  =\left(  -1\right)  ^{p+q}%
\det\left(  A_{\sim p,\sim q}\right)  .
\end{align*}
This proves Lemma \ref{lem.laplace.Apq}.
\end{proof}
\end{vershort}

\begin{verlong}
\begin{proof}
[Proof of Lemma \ref{lem.laplace.Apq}.]Let us use all notations introduced in
Lemma \ref{lem.laplace.gp}.

We have $p\in\left\{  1,2,\ldots,n\right\}  =\left[  n\right]  $. Hence,
$g_{p}$ is well-defined. Similarly, $g_{q}$ is well-defined. We have%
\[
\left(  g_{p}\left(  1\right)  ,g_{p}\left(  2\right)  ,\ldots,g_{p}\left(
n-1\right)  \right)  =\left(  1,2,\ldots,\widehat{p},\ldots,n\right)
\]
(by Lemma \ref{lem.laplace.gp} \textbf{(d)}) and
\[
\left(  g_{q}\left(  1\right)  ,g_{q}\left(  2\right)  ,\ldots,g_{q}\left(
n-1\right)  \right)  =\left(  1,2,\ldots,\widehat{q},\ldots,n\right)
\]
(by Lemma \ref{lem.laplace.gp} \textbf{(d)}, applied to $q$ instead of $p$).
Now, the definition of $A_{\sim p,\sim q}$ yields%
\begin{align}
A_{\sim p,\sim q}  &  =\operatorname*{sub}\nolimits_{1,2,\ldots,\widehat{p}%
,\ldots,n}^{1,2,\ldots,\widehat{q},\ldots,n}A=\operatorname*{sub}%
\nolimits_{1,2,\ldots,\widehat{p},\ldots,n}^{g_{q}\left(  1\right)
,g_{q}\left(  2\right)  ,\ldots,g_{q}\left(  n-1\right)  }A\nonumber\\
&  \ \ \ \ \ \ \ \ \ \ \left(  \text{since }\left(  1,2,\ldots,\widehat{q}%
,\ldots,n\right)  =\left(  g_{q}\left(  1\right)  ,g_{q}\left(  2\right)
,\ldots,g_{q}\left(  n-1\right)  \right)  \right) \nonumber\\
&  =\operatorname*{sub}\nolimits_{g_{p}\left(  1\right)  ,g_{p}\left(
2\right)  ,\ldots,g_{p}\left(  n-1\right)  }^{g_{q}\left(  1\right)
,g_{q}\left(  2\right)  ,\ldots,g_{q}\left(  n-1\right)  }A\nonumber\\
&  \ \ \ \ \ \ \ \ \ \ \left(  \text{since }\left(  1,2,\ldots,\widehat{p}%
,\ldots,n\right)  =\left(  g_{p}\left(  1\right)  ,g_{p}\left(  2\right)
,\ldots,g_{p}\left(  n-1\right)  \right)  \right) \nonumber\\
&  =\left(  a_{g_{p}\left(  x\right)  ,g_{q}\left(  y\right)  }\right)
_{1\leq x\leq n-1,\ 1\leq y\leq n-1}\nonumber\\
&  \ \ \ \ \ \ \ \ \ \ \left(  \text{by the definition of }\operatorname*{sub}%
\nolimits_{g_{p}\left(  1\right)  ,g_{p}\left(  2\right)  ,\ldots,g_{p}\left(
n-1\right)  }^{g_{q}\left(  1\right)  ,g_{q}\left(  2\right)  ,\ldots
,g_{q}\left(  n-1\right)  }A\right) \nonumber\\
&  =\left(  a_{g_{p}\left(  i\right)  ,g_{q}\left(  j\right)  }\right)
_{1\leq i\leq n-1,\ 1\leq j\leq n-1}\label{pf.lem.laplace.Apq.A}\\
&  \ \ \ \ \ \ \ \ \ \ \left(  \text{here, we renamed the index }\left(
x,y\right)  \text{ as }\left(  i,j\right)  \right)  .\nonumber
\end{align}

Also, $\left[  n\right]  $ is nonempty (since $p\in\left[  n\right]  $), and
thus we have $n>0$.

Now, let us recall the map $T:\left\{  \tau\in S_{n}\ \mid\ \tau\left(
n\right)  =n\right\}  \rightarrow\left\{  \tau\in S_{n}\ \mid\ \tau\left(
p\right)  =q\right\}  $ defined in Lemma \ref{lem.laplace.gp} \textbf{(g)}.
Lemma \ref{lem.laplace.gp} \textbf{(g)} says that this map $T$ is well-defined
and bijective. Every $\sigma\in\left\{  \tau\in S_{n}\ \mid\ \tau\left(
n\right)  =n\right\}  $ satisfies%
\begin{equation}
\left(  -1\right)  ^{T\left(  \sigma\right)  }=\left(  -1\right)  ^{p+q}%
\cdot\left(  -1\right)  ^{\sigma} \label{pf.lem.laplace.Apq.1}%
\end{equation}
\footnote{\textit{Proof of (\ref{pf.lem.laplace.Apq.1}):} Let $\sigma
\in\left\{  \tau\in S_{n}\ \mid\ \tau\left(  n\right)  =n\right\}  $. Applying
Lemma \ref{lem.laplace.gp} \textbf{(e)} to $q$ instead of $p$, we obtain%
\[
\left(  -1\right)  ^{g_{q}}=\left(  -1\right)  ^{n-q}=\left(  -1\right)
^{n+q}\ \ \ \ \ \ \ \ \ \ \left(  \text{since }n-q\equiv n+q\operatorname{mod}%
2\right)  .
\]
\par
The definition of $T\left(  \sigma\right)  $ yields $T\left(  \sigma\right)
=g_{q}\circ\sigma\circ\left(  g_{p}\right)  ^{-1}$. Thus,%
\[
\underbrace{T\left(  \sigma\right)  }_{=g_{q}\circ\sigma\circ\left(
g_{p}\right)  ^{-1}}\circ g_{p}=g_{q}\circ\sigma\circ\underbrace{\left(
g_{p}\right)  ^{-1}\circ g_{p}}_{=\operatorname*{id}}=g_{q}\circ\sigma,
\]
so that%
\begin{align*}
\left(  -1\right)  ^{T\left(  \sigma\right)  \circ g_{p}}  &  =\left(
-1\right)  ^{g_{q}\circ\sigma}=\underbrace{\left(  -1\right)  ^{g_{q}}%
}_{=\left(  -1\right)  ^{n+q}}\cdot\left(  -1\right)  ^{\sigma}%
\ \ \ \ \ \ \ \ \ \ \left(  \text{by (\ref{eq.sign.prod}), applied to }%
g_{q}\text{ and }\sigma\text{ instead of }\sigma\text{ and }\tau\right) \\
&  =\left(  -1\right)  ^{n+q}\cdot\left(  -1\right)  ^{\sigma}.
\end{align*}
Compared with%
\begin{align*}
\left(  -1\right)  ^{T\left(  \sigma\right)  \circ g_{p}}  &  =\left(
-1\right)  ^{T\left(  \sigma\right)  }\cdot\underbrace{\left(  -1\right)
^{g_{p}}}_{\substack{=\left(  -1\right)  ^{n-p}\\\text{(by Lemma
\ref{lem.laplace.gp} \textbf{(e)})}}}\ \ \ \ \ \ \ \ \ \ \left(  \text{by
(\ref{eq.sign.prod}), applied to }T\left(  \sigma\right)  \text{ and }%
g_{p}\text{ instead of }\sigma\text{ and }\tau\right) \\
&  =\left(  -1\right)  ^{T\left(  \sigma\right)  }\cdot\left(  -1\right)
^{n-p},
\end{align*}
this yields
\[
\left(  -1\right)  ^{T\left(  \sigma\right)  }\cdot\left(  -1\right)
^{n-p}=\left(  -1\right)  ^{n+q}\cdot\left(  -1\right)  ^{\sigma}.
\]
We can divide both sides of this equality by $\left(  -1\right)  ^{n-p}$
(since $\left(  -1\right)  ^{n-p}\in\left\{  1,-1\right\}  $ is clearly an
invertible integer), and thus we obtain%
\[
\left(  -1\right)  ^{T\left(  \sigma\right)  }=\dfrac{\left(  -1\right)
^{n+q}\cdot\left(  -1\right)  ^{\sigma}}{\left(  -1\right)  ^{n-p}%
}=\underbrace{\dfrac{\left(  -1\right)  ^{n+q}}{\left(  -1\right)  ^{n-p}}%
}_{\substack{=\left(  -1\right)  ^{\left(  n+q\right)  -\left(  n-p\right)
}=\left(  -1\right)  ^{p+q}\\\text{(since }\left(  n+q\right)  -\left(
n-p\right)  =p+q\text{)}}}\cdot\left(  -1\right)  ^{\sigma}=\left(  -1\right)
^{p+q}\cdot\left(  -1\right)  ^{\sigma}.
\]
This proves (\ref{pf.lem.laplace.Apq.1}).} and%
\begin{equation}
\prod_{\substack{i\in\left\{  1,2,\ldots,n\right\}  ;\\i\neq p}}a_{i,\left(
T\left(  \sigma\right)  \right)  \left(  i\right)  }=\prod_{i=1}^{n-1}%
a_{g_{p}\left(  i\right)  ,g_{q}\left(  \sigma\left(  i\right)  \right)  }
\label{pf.lem.laplace.Apq.2}%
\end{equation}
\footnote{\textit{Proof of (\ref{pf.lem.laplace.Apq.2}):} Let $\sigma
\in\left\{  \tau\in S_{n}\ \mid\ \tau\left(  n\right)  =n\right\}  $. Let us
recall the map $g_{p}^{\prime}:\left[  n-1\right]  \rightarrow\left[
n\right]  \setminus\left\{  p\right\}  $ introduced in Lemma
\ref{lem.laplace.gp} \textbf{(f)}. Lemma \ref{lem.laplace.gp} \textbf{(f)}
says that this map $g_{p}^{\prime}$ is well-defined and bijective. In other
words, $g_{p}^{\prime}$ is a bijection.
\par
Let $i\in\left[  n-1\right]  $. Then, $g_{p}^{\prime}\left(  i\right)
=g_{p}\left(  i\right)  $ (by the definition of $g_{p}^{\prime}$). Also, the
definition of $T$ yields $T\left(  \sigma\right)  =g_{q}\circ\sigma
\circ\left(  g_{p}\right)  ^{-1}$, so that%
\[
\left(  \underbrace{T\left(  \sigma\right)  }_{=g_{q}\circ\sigma\circ\left(
g_{p}\right)  ^{-1}}\right)  \left(  \underbrace{g_{p}^{\prime}\left(
i\right)  }_{=g_{p}\left(  i\right)  }\right)  =\left(  g_{q}\circ\sigma
\circ\left(  g_{p}\right)  ^{-1}\right)  \left(  g_{p}\left(  i\right)
\right)  =g_{q}\left(  \sigma\left(  \underbrace{\left(  g_{p}\right)
^{-1}\left(  g_{p}\left(  i\right)  \right)  }_{=i}\right)  \right)
=g_{q}\left(  \sigma\left(  i\right)  \right)  .
\]
\par
From $g_{p}^{\prime}\left(  i\right)  =g_{p}\left(  i\right)  $ and $\left(
T\left(  \sigma\right)  \right)  \left(  g_{p}^{\prime}\left(  i\right)
\right)  =g_{q}\left(  \sigma\left(  i\right)  \right)  $, we obtain%
\begin{equation}
a_{g_{p}^{\prime}\left(  i\right)  ,\left(  T\left(  \sigma\right)  \right)
\left(  g_{p}^{\prime}\left(  i\right)  \right)  }=a_{g_{p}\left(  i\right)
,g_{q}\left(  \sigma\left(  i\right)  \right)  }.
\label{pf.lem.laplace.Apq.2.pf.1}%
\end{equation}
\par
Now, let us forget that we fixed $i$. We thus have proven
(\ref{pf.lem.laplace.Apq.2.pf.1}) for every $i\in\left[  n-1\right]  $. But
now, we have%
\begin{align*}
&  \underbrace{\prod_{\substack{i\in\left\{  1,2,\ldots,n\right\}  ;\\i\neq
p}}}_{\substack{=\prod_{\substack{i\in\left[  n\right]  ;\\i\neq
p}}\\\text{(since }\left\{  1,2,\ldots,n\right\}  =\left[  n\right]  \text{)}%
}}a_{i,\left(  T\left(  \sigma\right)  \right)  \left(  i\right)  }\\
&  =\underbrace{\prod_{\substack{i\in\left[  n\right]  ;\\i\neq p}}}%
_{=\prod_{i\in\left[  n\right]  \setminus\left\{  p\right\}  }}a_{i,\left(
T\left(  \sigma\right)  \right)  \left(  i\right)  }=\prod_{i\in\left[
n\right]  \setminus\left\{  p\right\}  }a_{i,\left(  T\left(  \sigma\right)
\right)  \left(  i\right)  }=\underbrace{\prod_{i\in\left[  n-1\right]  }%
}_{\substack{=\prod_{i\in\left\{  1,2,\ldots,n-1\right\}  }\\\text{(since
}\left[  n-1\right]  =\left\{  1,2,\ldots,n-1\right\}  \text{)}}%
}\underbrace{a_{g_{p}^{\prime}\left(  i\right)  ,\left(  T\left(
\sigma\right)  \right)  \left(  g_{p}^{\prime}\left(  i\right)  \right)  }%
}_{\substack{=a_{g_{p}\left(  i\right)  ,g_{q}\left(  \sigma\left(  i\right)
\right)  }\\\text{(by (\ref{pf.lem.laplace.Apq.2.pf.1}))}}}\\
&  \ \ \ \ \ \ \ \ \ \ \left(
\begin{array}
[c]{c}%
\text{here, we have substituted }g_{p}^{\prime}\left(  i\right)  \text{ for
}i\text{, since}\\
g_{p}^{\prime}:\left[  n-1\right]  \rightarrow\left[  n\right]  \setminus
\left\{  p\right\}  \text{ is a bijection}%
\end{array}
\right) \\
&  =\underbrace{\prod_{i\in\left\{  1,2,\ldots,n-1\right\}  }}_{=\prod
_{i=1}^{n-1}}a_{g_{p}\left(  i\right)  ,g_{q}\left(  \sigma\left(  i\right)
\right)  }=\prod_{i=1}^{n-1}a_{g_{p}\left(  i\right)  ,g_{q}\left(
\sigma\left(  i\right)  \right)  }.
\end{align*}
This proves (\ref{pf.lem.laplace.Apq.2}).}.

Now,
\begin{align*}
&  \underbrace{\sum_{\substack{\sigma\in S_{n};\\\sigma\left(  p\right)  =q}%
}}_{=\sum_{\sigma\in\left\{  \tau\in S_{n}\ \mid\ \tau\left(  p\right)
=q\right\}  }}\left(  -1\right)  ^{\sigma}\prod_{\substack{i\in\left\{
1,2,\ldots,n\right\}  ;\\i\neq p}}a_{i,\sigma\left(  i\right)  }\\
&  =\sum_{\sigma\in\left\{  \tau\in S_{n}\ \mid\ \tau\left(  p\right)
=q\right\}  }\left(  -1\right)  ^{\sigma}\prod_{\substack{i\in\left\{
1,2,\ldots,n\right\}  ;\\i\neq p}}a_{i,\sigma\left(  i\right)  }\\
&  =\underbrace{\sum_{\sigma\in\left\{  \tau\in S_{n}\ \mid\ \tau\left(
n\right)  =n\right\}  }}_{=\sum_{\substack{\sigma\in S_{n};\\\sigma\left(
n\right)  =n}}}\underbrace{\left(  -1\right)  ^{T\left(  \sigma\right)  }%
}_{\substack{=\left(  -1\right)  ^{p+q}\cdot\left(  -1\right)  ^{\sigma
}\\\text{(by (\ref{pf.lem.laplace.Apq.1}))}}}\underbrace{\prod_{\substack{i\in
\left\{  1,2,\ldots,n\right\}  ;\\i\neq p}}a_{i,\left(  T\left(
\sigma\right)  \right)  \left(  i\right)  }}_{\substack{=\prod_{i=1}%
^{n-1}a_{g_{p}\left(  i\right)  ,g_{q}\left(  \sigma\left(  i\right)  \right)
}\\\text{(by (\ref{pf.lem.laplace.Apq.2}))}}}\\
&  \ \ \ \ \ \ \ \ \ \ \left(
\begin{array}
[c]{c}%
\text{here, we have substituted }T\left(  \sigma\right)  \text{ for }%
\sigma\text{ in the sum,}\\
\text{since the map }T:\left\{  \tau\in S_{n}\ \mid\ \tau\left(  n\right)
=n\right\}  \rightarrow\left\{  \tau\in S_{n}\ \mid\ \tau\left(  p\right)
=q\right\}
\end{array}
\right) \\
&  =\sum_{\substack{\sigma\in S_{n};\\\sigma\left(  n\right)  =n}}\left(
-1\right)  ^{p+q}\cdot\left(  -1\right)  ^{\sigma}\prod_{i=1}^{n-1}%
a_{g_{p}\left(  i\right)  ,g_{q}\left(  \sigma\left(  i\right)  \right)  }\\
&  =\left(  -1\right)  ^{p+q}\underbrace{\sum_{\substack{\sigma\in
S_{n};\\\sigma\left(  n\right)  =n}}\left(  -1\right)  ^{\sigma}\prod
_{i=1}^{n-1}a_{g_{p}\left(  i\right)  ,g_{q}\left(  \sigma\left(  i\right)
\right)  }}_{\substack{=\det\left(  \left(  a_{g_{p}\left(  i\right)
,g_{q}\left(  j\right)  }\right)  _{1\leq i\leq n-1,\ 1\leq j\leq n-1}\right)
\\\text{(by Lemma \ref{lem.laplace.lem}, applied to }a_{g_{p}\left(  i\right)
,g_{q}\left(  j\right)  }\text{ instead of }a_{i,j}\text{)}}}\\
&  =\left(  -1\right)  ^{p+q}\det\left(  \underbrace{\left(  a_{g_{p}\left(
i\right)  ,g_{q}\left(  j\right)  }\right)  _{1\leq i\leq n-1,\ 1\leq j\leq
n-1}}_{\substack{=A_{\sim p,\sim q}\\\text{(by (\ref{pf.lem.laplace.Apq.A}))}%
}}\right)  =\left(  -1\right)  ^{p+q}\det\left(  A_{\sim p,\sim q}\right)  .
\end{align*}
This proves Lemma \ref{lem.laplace.Apq}.
\end{proof}
\end{verlong}

Now, we can finally prove Theorem \ref{thm.laplace.gen}:

\begin{proof}
[Proof of Theorem \ref{thm.laplace.gen}.]\textbf{(a)} Let $p\in\left\{
1,2,\ldots,n\right\}  $. From (\ref{eq.det.eq.2}), we obtain%
\begin{align*}
\det A  &  =\sum_{\sigma\in S_{n}}\left(  -1\right)  ^{\sigma}\prod_{i=1}%
^{n}a_{i,\sigma\left(  i\right)  }\\
&  =\sum_{q\in\left\{  1,2,\ldots,n\right\}  }\sum_{\substack{\sigma\in
S_{n};\\\sigma\left(  p\right)  =q}}\left(  -1\right)  ^{\sigma}%
\underbrace{\prod_{i=1}^{n}}_{=\prod_{i\in\left\{  1,2,\ldots,n\right\}  }%
}a_{i,\sigma\left(  i\right)  }\\
&  \ \ \ \ \ \ \ \ \ \ \left(
\begin{array}
[c]{c}%
\text{because for every }\sigma\in S_{n}\text{, there exists}\\
\text{exactly one }q\in\left\{  1,2,\ldots,n\right\}  \text{ satisfying
}\sigma\left(  p\right)  =q
\end{array}
\right) \\
&  =\sum_{q\in\left\{  1,2,\ldots,n\right\}  }\sum_{\substack{\sigma\in
S_{n};\\\sigma\left(  p\right)  =q}}\left(  -1\right)  ^{\sigma}%
\underbrace{\prod_{i\in\left\{  1,2,\ldots,n\right\}  }a_{i,\sigma\left(
i\right)  }}_{\substack{=a_{p,\sigma\left(  p\right)  }\prod_{\substack{i\in
\left\{  1,2,\ldots,n\right\}  ;\\i\neq p}}a_{i,\sigma\left(  i\right)
}\\\text{(here, we have split off the factor for }i=p\text{ from the
product)}}}\\
&  =\underbrace{\sum_{q\in\left\{  1,2,\ldots,n\right\}  }}_{=\sum_{q=1}^{n}%
}\sum_{\substack{\sigma\in S_{n};\\\sigma\left(  p\right)  =q}}\left(
-1\right)  ^{\sigma}\underbrace{a_{p,\sigma\left(  p\right)  }}%
_{\substack{=a_{p,q}\\\text{(since }\sigma\left(  p\right)  =q\text{)}}%
}\prod_{\substack{i\in\left\{  1,2,\ldots,n\right\}  ;\\i\neq p}%
}a_{i,\sigma\left(  i\right)  }\\
&  =\sum_{q=1}^{n}\underbrace{\sum_{\substack{\sigma\in S_{n};\\\sigma\left(
p\right)  =q}}\left(  -1\right)  ^{\sigma}a_{p,q}\prod_{\substack{i\in\left\{
1,2,\ldots,n\right\}  ;\\i\neq p}}a_{i,\sigma\left(  i\right)  }}%
_{=a_{p,q}\sum_{\substack{\sigma\in S_{n};\\\sigma\left(  p\right)
=q}}\left(  -1\right)  ^{\sigma}\prod_{\substack{i\in\left\{  1,2,\ldots
,n\right\}  ;\\i\neq p}}a_{i,\sigma\left(  i\right)  }}\\
&  =\sum_{q=1}^{n}a_{p,q}\underbrace{\sum_{\substack{\sigma\in S_{n}%
;\\\sigma\left(  p\right)  =q}}\left(  -1\right)  ^{\sigma}\prod
_{\substack{i\in\left\{  1,2,\ldots,n\right\}  ;\\i\neq p}}a_{i,\sigma\left(
i\right)  }}_{\substack{=\left(  -1\right)  ^{p+q}\det\left(  A_{\sim p,\sim
q}\right)  \\\text{(by Lemma \ref{lem.laplace.Apq})}}}\\
&  =\sum_{q=1}^{n}\underbrace{a_{p,q}\left(  -1\right)  ^{p+q}}_{=\left(
-1\right)  ^{p+q}a_{p,q}}\det\left(  A_{\sim p,\sim q}\right)  =\sum_{q=1}%
^{n}\left(  -1\right)  ^{p+q}a_{p,q}\det\left(  A_{\sim p,\sim q}\right)  .
\end{align*}
This proves Theorem \ref{thm.laplace.gen} \textbf{(a)}.

\textbf{(b)} Let $q\in\left\{  1,2,\ldots,n\right\}  $. From
(\ref{eq.det.eq.2}), we obtain%
\begin{align*}
\det A  &  =\sum_{\sigma\in S_{n}}\left(  -1\right)  ^{\sigma}\prod_{i=1}%
^{n}a_{i,\sigma\left(  i\right)  }\\
&  =\sum_{p\in\left\{  1,2,\ldots,n\right\}  }\sum_{\substack{\sigma\in
S_{n};\\\sigma^{-1}\left(  q\right)  =p}}\left(  -1\right)  ^{\sigma
}\underbrace{\prod_{i=1}^{n}}_{=\prod_{i\in\left\{  1,2,\ldots,n\right\}  }%
}a_{i,\sigma\left(  i\right)  }\\
&  \ \ \ \ \ \ \ \ \ \ \left(
\begin{array}
[c]{c}%
\text{because for every }\sigma\in S_{n}\text{, there exists}\\
\text{exactly one }p\in\left\{  1,2,\ldots,n\right\}  \text{ satisfying
}\sigma^{-1}\left(  q\right)  =p
\end{array}
\right) \\
&  =\sum_{p\in\left\{  1,2,\ldots,n\right\}  }\underbrace{\sum
_{\substack{\sigma\in S_{n};\\\sigma^{-1}\left(  q\right)  =p}}}%
_{\substack{=\sum_{\substack{\sigma\in S_{n};\\\sigma\left(  p\right)
=q}}\\\text{(because for any }\sigma\in S_{n}\text{,}\\\text{the statement
}\left(  \sigma^{-1}\left(  q\right)  =p\right)  \\\text{is equivalent to
the}\\\text{statement }\left(  \sigma\left(  p\right)  =q\right)  \text{)}%
}}\left(  -1\right)  ^{\sigma}\underbrace{\prod_{i\in\left\{  1,2,\ldots
,n\right\}  }a_{i,\sigma\left(  i\right)  }}_{\substack{=a_{p,\sigma\left(
p\right)  }\prod_{\substack{i\in\left\{  1,2,\ldots,n\right\}  ;\\i\neq
p}}a_{i,\sigma\left(  i\right)  }\\\text{(here, we have split off
the}\\\text{factor for }i=p\text{ from the product)}}}\\
&  =\underbrace{\sum_{p\in\left\{  1,2,\ldots,n\right\}  }}_{=\sum_{p=1}^{n}%
}\sum_{\substack{\sigma\in S_{n};\\\sigma\left(  p\right)  =q}}\left(
-1\right)  ^{\sigma}\underbrace{a_{p,\sigma\left(  p\right)  }}%
_{\substack{=a_{p,q}\\\text{(since }\sigma\left(  p\right)  =q\text{)}}%
}\prod_{\substack{i\in\left\{  1,2,\ldots,n\right\}  ;\\i\neq p}%
}a_{i,\sigma\left(  i\right)  }\\
&  =\sum_{p=1}^{n}\underbrace{\sum_{\substack{\sigma\in S_{n};\\\sigma\left(
p\right)  =q}}\left(  -1\right)  ^{\sigma}a_{p,q}\prod_{\substack{i\in\left\{
1,2,\ldots,n\right\}  ;\\i\neq p}}a_{i,\sigma\left(  i\right)  }}%
_{=a_{p,q}\sum_{\substack{\sigma\in S_{n};\\\sigma\left(  p\right)
=q}}\left(  -1\right)  ^{\sigma}\prod_{\substack{i\in\left\{  1,2,\ldots
,n\right\}  ;\\i\neq p}}a_{i,\sigma\left(  i\right)  }}\\
&  =\sum_{p=1}^{n}a_{p,q}\underbrace{\sum_{\substack{\sigma\in S_{n}%
;\\\sigma\left(  p\right)  =q}}\left(  -1\right)  ^{\sigma}\prod
_{\substack{i\in\left\{  1,2,\ldots,n\right\}  ;\\i\neq p}}a_{i,\sigma\left(
i\right)  }}_{\substack{=\left(  -1\right)  ^{p+q}\det\left(  A_{\sim p,\sim
q}\right)  \\\text{(by Lemma \ref{lem.laplace.Apq})}}}\\
&  =\sum_{p=1}^{n}\underbrace{a_{p,q}\left(  -1\right)  ^{p+q}}_{=\left(
-1\right)  ^{p+q}a_{p,q}}\det\left(  A_{\sim p,\sim q}\right)  =\sum_{p=1}%
^{n}\left(  -1\right)  ^{p+q}a_{p,q}\det\left(  A_{\sim p,\sim q}\right)  .
\end{align*}
This proves Theorem \ref{thm.laplace.gen} \textbf{(b)}.
\end{proof}

Let me make three simple observations (which can easily be checked by the reader):

\begin{itemize}
\item Theorem \ref{thm.laplace.gen} \textbf{(b)} could be (alternatively)
proven using Theorem \ref{thm.laplace.gen} \textbf{(a)} (applied to $A^{T}$
and $a_{q,p}$ instead of $A$ and $a_{p,q}$) and Exercise \ref{exe.ps4.4}.

\item Theorem \ref{thm.laplace.pre} is a particular case of Theorem
\ref{thm.laplace.gen} \textbf{(a)}.

\item Corollary \ref{cor.laplace.pre.col} is a particular case of Theorem
\ref{thm.laplace.gen} \textbf{(b)}.
\end{itemize}

\begin{remark}
Some books use Laplace expansion to define the notion of a determinant. For
example, one can define the determinant of a square matrix recursively, by
setting the determinant of the $0\times0$-matrix to be $1$, and defining the
determinant of an $n\times n$-matrix $A=\left(  a_{i,j}\right)  _{1\leq i\leq
n,\ 1\leq j\leq n}$ (with $n>0$) to be $\sum_{q=1}^{n}\left(  -1\right)
^{1+q}a_{1,q}\det\left(  A_{\sim1,\sim q}\right)  $ (assuming that
determinants of $\left(  n-1\right)  \times\left(  n-1\right)  $-matrices such
as $A_{\sim1,\sim q}$ are already defined). Of course, this leads to the same
notion of determinant as the one we are using, because of Theorem
\ref{thm.laplace.gen} \textbf{(a)}.
\end{remark}

\subsection{\label{sect.tridiag}Tridiagonal determinants}

In this section, we shall study the so-called \textit{tridiagonal matrices}: a
class of matrices whose all entries are zero everywhere except in the
\textquotedblleft direct proximity\textquotedblright\ of the diagonal (more
specifically: on the diagonal and \textquotedblleft one level below and one
level above\textquotedblright). We shall find recursive formulas for the
determinants of these matrices. These formulas are a simple example of an
application of Laplace expansion, but also interesting in their own right.

\begin{definition}
\label{def.tridiag}Let $n\in\mathbb{N}$. Let $a_{1},a_{2},\ldots,a_{n}$ be $n$
elements of $\mathbb{K}$. Let $b_{1},b_{2},\ldots,b_{n-1}$ be $n-1$ elements
of $\mathbb{K}$ (where we take the position that \textquotedblleft$-1$
elements of $\mathbb{K}$\textquotedblright\ means \textquotedblleft no
elements of $\mathbb{K}$\textquotedblright). Let $c_{1},c_{2},\ldots,c_{n-1}$
be $n-1$ elements of $\mathbb{K}$. We now set%
\[
A=\left(
\begin{array}
[c]{ccccccc}%
a_{1} & b_{1} & 0 & \cdots & 0 & 0 & 0\\
c_{1} & a_{2} & b_{2} & \cdots & 0 & 0 & 0\\
0 & c_{2} & a_{3} & \cdots & 0 & 0 & 0\\
\vdots & \vdots & \vdots & \ddots & \vdots & \vdots & \vdots\\
0 & 0 & 0 & \cdots & a_{n-2} & b_{n-2} & 0\\
0 & 0 & 0 & \cdots & c_{n-2} & a_{n-1} & b_{n-1}\\
0 & 0 & 0 & \cdots & 0 & c_{n-1} & a_{n}%
\end{array}
\right)  .
\]
(More formally,%
\[
A=\left(
\begin{cases}
a_{i}, & \text{if }i=j;\\
b_{i}, & \text{if }i=j-1;\\
c_{j}, & \text{if }i=j+1;\\
0, & \text{otherwise}%
\end{cases}
\right)  _{1\leq i\leq n,\ 1\leq j\leq n}.
\]
)

The matrix $A$ is called a \textit{tridiagonal matrix}.

We shall keep the notations $n$, $a_{1},a_{2},\ldots,a_{n}$, $b_{1}%
,b_{2},\ldots,b_{n-1}$, $c_{1},c_{2},\ldots,c_{n-1}$ and $A$ fixed for the
rest of Section \ref{sect.tridiag}.
\end{definition}

Playing around with small examples, one soon notices that the determinants of
tridiagonal matrices are too complicated to have neat explicit formulas in
full generality. For $n\in\left\{  0,1,2,3\right\}  $, the determinants look
as follows:%
\begin{align*}
\det A  &  =\det\left(  \text{the }0\times0\text{-matrix}\right)
=1\ \ \ \ \ \ \ \ \ \ \text{if }n=0;\\
\det A  &  =\det\left(
\begin{array}
[c]{c}%
a_{1}%
\end{array}
\right)  =a_{1}\ \ \ \ \ \ \ \ \ \ \text{if }n=1;\\
\det A  &  =\det\left(
\begin{array}
[c]{cc}%
a_{1} & b_{1}\\
c_{1} & a_{2}%
\end{array}
\right)  =a_{1}a_{2}-b_{1}c_{1}\ \ \ \ \ \ \ \ \ \ \text{if }n=2;\\
\det A  &  =\det\left(
\begin{array}
[c]{ccc}%
a_{1} & b_{1} & 0\\
c_{1} & a_{2} & b_{2}\\
0 & c_{2} & a_{3}%
\end{array}
\right)  =a_{1}a_{2}a_{3}-a_{1}b_{2}c_{2}-a_{3}b_{1}c_{1}%
\ \ \ \ \ \ \ \ \ \ \text{if }n=3.
\end{align*}
(And these formulas get more complicated the larger $n$ becomes.) However, the
many zeroes present in a tridiagonal matrix make it easy to find a recursive
formula for its determinant using Laplace expansion:

\begin{proposition}
\label{prop.tridiag.rec}For every two elements $x$ and $y$ of $\left\{
0,1,\ldots,n\right\}  $ satisfying $x\leq y$, we let $A_{x,y}$ be the $\left(
y-x\right)  \times\left(  y-x\right)  $-matrix%
\[
\left(
\begin{array}
[c]{ccccccc}%
a_{x+1} & b_{x+1} & 0 & \cdots & 0 & 0 & 0\\
c_{x+1} & a_{x+2} & b_{x+2} & \cdots & 0 & 0 & 0\\
0 & c_{x+2} & a_{x+3} & \cdots & 0 & 0 & 0\\
\vdots & \vdots & \vdots & \ddots & \vdots & \vdots & \vdots\\
0 & 0 & 0 & \cdots & a_{y-2} & b_{y-2} & 0\\
0 & 0 & 0 & \cdots & c_{y-2} & a_{y-1} & b_{y-1}\\
0 & 0 & 0 & \cdots & 0 & c_{y-1} & a_{y}%
\end{array}
\right)  =\operatorname*{sub}\nolimits_{x+1,x+2,\ldots,y}^{x+1,x+2,\ldots
,y}A.
\]

\textbf{(a)} We have $\det\left(  A_{x,x}\right)  =1$ for every $x\in\left\{
0,1,\ldots,n\right\}  $.

\textbf{(b)} We have $\det\left(  A_{x,x+1}\right)  =a_{x+1}$ for every
$x\in\left\{  0,1,\ldots,n-1\right\}  $.

\textbf{(c)} For every $x\in\left\{  0,1,\ldots,n\right\}  $ and $y\in\left\{
0,1,\ldots,n\right\}  $ satisfying $x\leq y-2$, we have
\[
\det\left(  A_{x,y}\right)  =a_{y}\det\left(  A_{x,y-1}\right)  -b_{y-1}%
c_{y-1}\det\left(  A_{x,y-2}\right)  .
\]

\textbf{(d)} For every $x\in\left\{  0,1,\ldots,n\right\}  $ and $y\in\left\{
0,1,\ldots,n\right\}  $ satisfying $x\leq y-2$, we have
\[
\det\left(  A_{x,y}\right)  =a_{x+1}\det\left(  A_{x+1,y}\right)
-b_{x+1}c_{x+1}\det\left(  A_{x+2,y}\right)  .
\]

\textbf{(e)} We have $A=A_{0,n}$.
\end{proposition}

\begin{proof}
[Proof of Proposition \ref{prop.tridiag.rec}.]\textbf{(e)} The definition of
$A_{0,n}$ yields
\[
A_{0,n}=\left(
\begin{array}
[c]{ccccccc}%
a_{1} & b_{1} & 0 & \cdots & 0 & 0 & 0\\
c_{1} & a_{2} & b_{2} & \cdots & 0 & 0 & 0\\
0 & c_{2} & a_{3} & \cdots & 0 & 0 & 0\\
\vdots & \vdots & \vdots & \ddots & \vdots & \vdots & \vdots\\
0 & 0 & 0 & \cdots & a_{n-2} & b_{n-2} & 0\\
0 & 0 & 0 & \cdots & c_{n-2} & a_{n-1} & b_{n-1}\\
0 & 0 & 0 & \cdots & 0 & c_{n-1} & a_{n}%
\end{array}
\right)  =A.
\]
This proves Proposition \ref{prop.tridiag.rec} \textbf{(e)}.

\textbf{(a)} Let $x\in\left\{  0,1,\ldots,n\right\}  $. Then, $A_{x,x}$ is an
$\left(  x-x\right)  \times\left(  x-x\right)  $-matrix, thus a $0\times
0$-matrix. Hence, its determinant is $\det\left(  A_{x,x}\right)  =1$. This
proves Proposition \ref{prop.tridiag.rec} \textbf{(a)}.

\textbf{(b)} Let $x\in\left\{  0,1,\ldots,n-1\right\}  $. The definition of
$A_{x,x+1}$ shows that $A_{x,x+1}$ is the $1\times1$-matrix $\left(
\begin{array}
[c]{c}%
a_{x+1}%
\end{array}
\right)  $. Hence, $\det\left(  A_{x,x+1}\right)  =\det\left(
\begin{array}
[c]{c}%
a_{x+1}%
\end{array}
\right)  =a_{x+1}$. This proves Proposition \ref{prop.tridiag.rec}
\textbf{(b)}.

\textbf{(c)} Let $x\in\left\{  0,1,\ldots,n\right\}  $ and $y\in\left\{
0,1,\ldots,n\right\}  $ be such that $x\leq y-2$. We have%
\begin{equation}
A_{x,y}=\left(
\begin{array}
[c]{ccccccc}%
a_{x+1} & b_{x+1} & 0 & \cdots & 0 & 0 & 0\\
c_{x+1} & a_{x+2} & b_{x+2} & \cdots & 0 & 0 & 0\\
0 & c_{x+2} & a_{x+3} & \cdots & 0 & 0 & 0\\
\vdots & \vdots & \vdots & \ddots & \vdots & \vdots & \vdots\\
0 & 0 & 0 & \cdots & a_{y-2} & b_{y-2} & 0\\
0 & 0 & 0 & \cdots & c_{y-2} & a_{y-1} & b_{y-1}\\
0 & 0 & 0 & \cdots & 0 & c_{y-1} & a_{y}%
\end{array}
\right)  . \label{pf.prop.tridiag.rec.c.Axy}%
\end{equation}
This is a $\left(  y-x\right)  \times\left(  y-x\right)  $-matrix. If we cross
out its $\left(  y-x\right)  $-th row (i.e., its last row) and its $\left(
y-x\right)  $-th column (i.e., its last column), then we obtain $A_{x,y-1}$.
In other words, $\left(  A_{x,y}\right)  _{\sim\left(  y-x\right)
,\sim\left(  y-x\right)  }=A_{x,y-1}$.

Let us write the matrix $A_{x,y}$ in the form $A_{x,y}=\left(  u_{i,j}\right)
_{1\leq i\leq y-x,\ 1\leq j\leq y-x}$. Thus,
\begin{align*}
&  \left(  u_{y-x,1},u_{y-x,2},\ldots,u_{y-x,y-x}\right) \\
&  =\left(  \text{the last row of the matrix }A_{x,y}\right)  =\left(
0,0,\ldots,0,c_{y-1},a_{y}\right)  .
\end{align*}
In other words, we have%
\begin{align}
&  \left(  u_{y-x,q}=0\ \ \ \ \ \ \ \ \ \ \text{for every }q\in\left\{
1,2,\ldots,y-x-2\right\}  \right)  ,\label{pf.prop.tridiag.rec.c.1}\\
&  u_{y-x,y-x-1}=c_{y-1},\ \ \ \ \ \ \ \ \ \ \text{and}\nonumber\\
&  u_{y-x,y-x}=a_{y}.\nonumber
\end{align}

Now, Laplace expansion along the $\left(  y-x\right)  $-th row (or, more
precisely, Theorem \ref{thm.laplace.gen} \textbf{(a)}, applied to $y-x$,
$A_{x,y}$, $u_{i,j}$ and $y-x$ instead of $n$, $A$, $a_{i,j}$ and $p$) yields%
\begin{align}
\det\left(  A_{x,y}\right)   &  =\sum_{q=1}^{y-x}\left(  -1\right)  ^{\left(
y-x\right)  +q}u_{y-x,q}\det\left(  \left(  A_{x,y}\right)  _{\sim\left(
y-x\right)  ,\sim q}\right) \nonumber\\
&  =\sum_{q=1}^{y-x-2}\left(  -1\right)  ^{\left(  y-x\right)  +q}%
\underbrace{u_{y-x,q}}_{\substack{=0\\\text{(by (\ref{pf.prop.tridiag.rec.c.1}%
))}}}\det\left(  \left(  A_{x,y}\right)  _{\sim\left(  y-x\right)  ,\sim
q}\right) \nonumber\\
&  \ \ \ \ \ \ \ \ \ \ +\underbrace{\left(  -1\right)  ^{\left(  y-x\right)
+\left(  y-x-1\right)  }}_{=-1}\underbrace{u_{y-x,y-x-1}}_{=c_{y-1}}%
\det\left(  \left(  A_{x,y}\right)  _{\sim\left(  y-x\right)  ,\sim\left(
y-x-1\right)  }\right) \nonumber\\
&  \ \ \ \ \ \ \ \ \ \ +\underbrace{\left(  -1\right)  ^{\left(  y-x\right)
+\left(  y-x\right)  }}_{=1}\underbrace{u_{y-x,y-x}}_{=a_{y}}\det\left(
\underbrace{\left(  A_{x,y}\right)  _{\sim\left(  y-x\right)  ,\sim\left(
y-x\right)  }}_{=A_{x,y-1}}\right) \nonumber\\
&  \ \ \ \ \ \ \ \ \ \ \left(  \text{since }y-x\geq2\text{ (since }x\leq
y-2\text{)}\right) \nonumber\\
&  =\underbrace{\sum_{q=1}^{y-x-2}\left(  -1\right)  ^{\left(  y-x\right)
+q}0\det\left(  \left(  A_{x,y}\right)  _{\sim\left(  y-x\right)  ,\sim
q}\right)  }_{=0}\nonumber\\
&  \ \ \ \ \ \ \ \ \ \ -c_{y-1}\det\left(  \left(  A_{x,y}\right)
_{\sim\left(  y-x\right)  ,\sim\left(  y-x-1\right)  }\right)  +a_{y}%
\det\left(  A_{x,y-1}\right) \nonumber\\
&  =-c_{y-1}\det\left(  \left(  A_{x,y}\right)  _{\sim\left(  y-x\right)
,\sim\left(  y-x-1\right)  }\right)  +a_{y}\det\left(  A_{x,y-1}\right)  .
\label{pf.prop.tridiag.rec.c.3}%
\end{align}

Now, let $B=\left(  A_{x,y}\right)  _{\sim\left(  y-x\right)  ,\sim\left(
y-x-1\right)  }$. Thus, (\ref{pf.prop.tridiag.rec.c.3}) becomes%
\begin{align}
\det\left(  A_{x,y}\right)   &  =-c_{y-1}\det\left(  \underbrace{\left(
A_{x,y}\right)  _{\sim\left(  y-x\right)  ,\sim\left(  y-x-1\right)  }}%
_{=B}\right)  +a_{y}\det\left(  A_{x,y-1}\right) \nonumber\\
&  =-c_{y-1}\det B+a_{y}\det\left(  A_{x,y-1}\right)  .
\label{pf.prop.tridiag.rec.c.3a}%
\end{align}

Now,%
\begin{align}
B  &  =\left(  A_{x,y}\right)  _{\sim\left(  y-x\right)  ,\sim\left(
y-x-1\right)  }\nonumber\\
&  =\left(
\begin{array}
[c]{ccccccc}%
a_{x+1} & b_{x+1} & 0 & \cdots & 0 & 0 & 0\\
c_{x+1} & a_{x+2} & b_{x+2} & \cdots & 0 & 0 & 0\\
0 & c_{x+2} & a_{x+3} & \cdots & 0 & 0 & 0\\
\vdots & \vdots & \vdots & \ddots & \vdots & \vdots & \vdots\\
0 & 0 & 0 & \cdots & a_{y-3} & b_{y-3} & 0\\
0 & 0 & 0 & \cdots & c_{y-3} & a_{y-2} & 0\\
0 & 0 & 0 & \cdots & 0 & c_{y-2} & b_{y-1}%
\end{array}
\right)  \ \ \ \ \ \ \ \ \ \ \left(  \text{because of
(\ref{pf.prop.tridiag.rec.c.Axy})}\right)  . \label{pf.prop.tridiag.rec.c.5}%
\end{align}

Now, let us write the matrix $B$ in the form $B=\left(  v_{i,j}\right)
_{1\leq i\leq y-x-1,\ 1\leq j\leq y-x-1}$. Thus,
\begin{align*}
&  \left(  v_{1,y-x-1},v_{2,y-x-1},\ldots,v_{y-x-1,y-x-1}\right)  ^{T}\\
&  =\left(  \text{the last column of the matrix }B\right)  =\left(
0,0,\ldots,0,b_{y-1}\right)  ^{T}%
\end{align*}
(because of (\ref{pf.prop.tridiag.rec.c.5})). In other words, we have%
\begin{align}
&  \left(  v_{p,y-x-1}=0\ \ \ \ \ \ \ \ \ \ \text{for every }p\in\left\{
1,2,\ldots,y-x-2\right\}  \right)  ,\ \ \ \ \ \ \ \ \ \ \text{and}%
\label{pf.prop.tridiag.rec.c.7}\\
&  v_{y-x-1,y-x-1}=b_{y-1}.\nonumber
\end{align}
Now, Laplace expansion along the $\left(  y-x-1\right)  $-th column (or, more
precisely, Theorem \ref{thm.laplace.gen} \textbf{(b)}, applied to $y-x-1$,
$B$, $v_{i,j}$ and $y-x-1$ instead of $n$, $A$, $a_{i,j}$ and $q$) yields%
\begin{align}
\det B  &  =\sum_{p=1}^{y-x-1}\left(  -1\right)  ^{p+\left(  y-x-1\right)
}v_{p,y-x-1}\det\left(  B_{\sim p,\sim\left(  y-x-1\right)  }\right)
\nonumber\\
&  =\sum_{p=1}^{y-x-2}\left(  -1\right)  ^{p+\left(  y-x-1\right)
}\underbrace{v_{p,y-x-1}}_{\substack{=0\\\text{(by
(\ref{pf.prop.tridiag.rec.c.7}))}}}\det\left(  B_{\sim p,\sim\left(
y-x-1\right)  }\right) \nonumber\\
&  \ \ \ \ \ \ \ \ \ \ +\underbrace{\left(  -1\right)  ^{\left(  y-x-1\right)
+\left(  y-x-1\right)  }}_{=1}\underbrace{v_{y-x-1,y-x-1}}_{=b_{y-1}}%
\det\left(  B_{\sim\left(  y-x-1\right)  ,\sim\left(  y-x-1\right)  }\right)
\nonumber\\
&  \ \ \ \ \ \ \ \ \ \ \left(  \text{since }y-x-1\geq1\text{ (since }x\leq
y-2\text{)}\right) \nonumber\\
&  =\underbrace{\sum_{p=1}^{y-x-2}\left(  -1\right)  ^{p+\left(  y-x-1\right)
}0\det\left(  B_{\sim p,\sim\left(  y-x-1\right)  }\right)  }_{=0}+b_{y-1}%
\det\left(  B_{\sim\left(  y-x-1\right)  ,\sim\left(  y-x-1\right)  }\right)
\nonumber\\
&  =b_{y-1}\det\left(  B_{\sim\left(  y-x-1\right)  ,\sim\left(  y-x-1\right)
}\right)  . \label{pf.prop.tridiag.rec.c.9}%
\end{align}
Finally, a look at (\ref{pf.prop.tridiag.rec.c.5}) reveals that%
\[
B_{\sim\left(  y-x-1\right)  ,\sim\left(  y-x-1\right)  }=\left(
\begin{array}
[c]{ccccccc}%
a_{x+1} & b_{x+1} & 0 & \cdots & 0 & 0 & 0\\
c_{x+1} & a_{x+2} & b_{x+2} & \cdots & 0 & 0 & 0\\
0 & c_{x+2} & a_{x+3} & \cdots & 0 & 0 & 0\\
\vdots & \vdots & \vdots & \ddots & \vdots & \vdots & \vdots\\
0 & 0 & 0 & \cdots & a_{y-4} & b_{y-4} & 0\\
0 & 0 & 0 & \cdots & c_{y-4} & a_{y-3} & b_{y-3}\\
0 & 0 & 0 & \cdots & 0 & c_{y-3} & a_{y-2}%
\end{array}
\right)  =A_{x,y-2}.
\]
Hence, (\ref{pf.prop.tridiag.rec.c.9}) becomes%
\[
\det B=b_{y-1}\det\left(  \underbrace{B_{\sim\left(  y-x-1\right)
,\sim\left(  y-x-1\right)  }}_{=A_{x,y-2}}\right)  =b_{y-1}\det\left(
A_{x,y-2}\right)  .
\]
Therefore, (\ref{pf.prop.tridiag.rec.c.3a}) becomes%
\begin{align*}
\det\left(  A_{x,y}\right)   &  =-c_{y-1}\underbrace{\det B}_{=b_{y-1}%
\det\left(  A_{x,y-2}\right)  }+a_{y}\det\left(  A_{x,y-1}\right) \\
&  =-c_{y-1}b_{y-1}\det\left(  A_{x,y-2}\right)  +a_{y}\det\left(
A_{x,y-1}\right) \\
&  =a_{y}\det\left(  A_{x,y-1}\right)  -b_{y-1}c_{y-1}\det\left(
A_{x,y-2}\right)  .
\end{align*}
This proves Proposition \ref{prop.tridiag.rec} \textbf{(c)}.

\textbf{(d)} The proof of Proposition \ref{prop.tridiag.rec} \textbf{(d)} is
similar to the proof of Proposition \ref{prop.tridiag.rec} \textbf{(c)}. The
main difference is that we now have to perform Laplace expansion along the
$1$-st row (instead of the $\left(  y-x\right)  $-th row) and then Laplace
expansion along the $1$-st column (instead of the $\left(  y-x-1\right)  $-th column).
\end{proof}

Proposition \ref{prop.tridiag.rec} gives us two fast recursive algorithms to
compute $\det A$:

The first algorithm proceeds by recursively computing $\det\left(
A_{0,m}\right)  $ for every $m\in\left\{  0,1,\ldots,n\right\}  $. This is
done using Proposition \ref{prop.tridiag.rec} \textbf{(a)} (for $m=0$),
Proposition \ref{prop.tridiag.rec} \textbf{(b)} (for $m=1$) and Proposition
\ref{prop.tridiag.rec} \textbf{(c)} (to find $\det\left(  A_{0,m}\right)  $
for $m\geq2$ in terms of $\det\left(  A_{0,m-1}\right)  $ and $\det\left(
A_{0,m-2}\right)  $). The final value $\det\left(  A_{0,n}\right)  $ is $\det
A$ (by Proposition \ref{prop.tridiag.rec} \textbf{(e)}).

The second algorithm proceeds by recursively computing $\det\left(
A_{m,n}\right)  $ for every $m\in\left\{  0,1,\ldots,n\right\}  $. This
recursion goes backwards: We start with $m=n$ (where we use Proposition
\ref{prop.tridiag.rec} \textbf{(a)}), then turn to $m=n-1$ (using Proposition
\ref{prop.tridiag.rec} \textbf{(b)}), and then go further and further down
(using Proposition \ref{prop.tridiag.rec} \textbf{(d)} to compute $\det\left(
A_{m,n}\right)  $ in terms of $\det\left(  A_{m+1,n}\right)  $ and
$\det\left(  A_{m+2,n}\right)  $).

So we have two different recursive algorithms leading to one and the same
result. Whenever you have such a thing, you can package up the equivalence of
the two algorithms as an exercise, and try to make it less easy by covering up
the actual goal of the algorithms (in our case, computing $\det A$). In our
case, this leads to the following exercise:

\begin{exercise}
\label{exe.tridiag.isl}Let $n\in\mathbb{N}$. Let $a_{1},a_{2},\ldots,a_{n}$ be
$n$ elements of $\mathbb{K}$. Let $b_{1},b_{2},\ldots,b_{n-1}$ be $n-1$
elements of $\mathbb{K}$.

Define a sequence $\left(  u_{0},u_{1},\ldots,u_{n}\right)  $ of elements of
$\mathbb{K}$ recursively by setting $u_{0}=1$, $u_{1}=a_{1}$ and%
\[
u_{i}=a_{i}u_{i-1}-b_{i-1}u_{i-2}\ \ \ \ \ \ \ \ \ \ \text{for every }%
i\in\left\{  2,3,\ldots,n\right\}  .
\]

Define a sequence $\left(  v_{0},v_{1},\ldots,v_{n}\right)  $ of elements of
$\mathbb{K}$ recursively by setting $v_{0}=1$, $v_{1}=a_{n}$ and%
\[
v_{i}=a_{n-i+1}v_{i-1}-b_{n-i+1}v_{i-2}\ \ \ \ \ \ \ \ \ \ \text{for every
}i\in\left\{  2,3,\ldots,n\right\}  .
\]

Prove that $u_{n}=v_{n}$.
\end{exercise}

This exercise generalizes
\href{http://www.artofproblemsolving.com/community/c6h597118p3543340}{IMO
Shortlist 2013 problem A1}\footnote{I have a suspicion that
\href{http://www.artofproblemsolving.com/community/c6h355915}{IMO Shortlist
2009 problem C3} also can be viewed as an equality between two recursive ways
to compute a determinant; but this determinant seems to be harder to find (I
don't think it can be obtained from Proposition \ref{prop.tridiag.rec}).}.

Our recursive algorithms for computing $\det A$ also yield another
observation: The determinant $\det A$ depends not on the $2\left(  n-1\right)
$ elements \newline$b_{1},b_{2},\ldots,b_{n-1},c_{1},c_{2},\ldots,c_{n-1}$ but
only on the products $b_{1}c_{1},b_{2}c_{2},\ldots,b_{n-1}c_{n-1}$.

\begin{exercise}
\label{exe.tridiag.cf}Define $A_{x,y}$ as in Proposition
\ref{prop.tridiag.rec}. Prove that%
\[
\dfrac{\det A}{\det\left(  A_{1,n}\right)  }=a_{1}-\dfrac{b_{1}c_{1}}%
{a_{2}-\dfrac{b_{2}c_{2}}{a_{3}-\dfrac{b_{3}c_{3}}{%
\begin{array}
[c]{ccc}%
a_{4}- &  & \\
& \ddots & \\
&  & -\dfrac{b_{n-2}c_{n-2}}{a_{n-1}-\dfrac{b_{n-1}c_{n-1}}{a_{n}}}%
\end{array}
}}},
\]
provided that all denominators in this equality are invertible.
\end{exercise}

\begin{exercise}
\label{exe.tridiag.fib}Assume that $a_{i}=1$ for all $i\in\left\{
1,2,\ldots,n\right\}  $. Also, assume that $b_{i}=1$ and $c_{i}=-1$ for all
$i\in\left\{  1,2,\ldots,n-1\right\}  $. Let $\left(  f_{0},f_{1},f_{2}%
,\ldots\right)  $ be the Fibonacci sequence (defined as in Chapter
\ref{chp.recur}). Show that $\det A=f_{n+1}$.
\end{exercise}

\begin{remark}
\label{rmk.tridiag.fib-cont}Consider once again the Fibonacci sequence
$\left(  f_{0},f_{1},f_{2},\ldots\right)  $ (defined as in Chapter
\ref{chp.recur}). Let $n$ be a positive integer. Combining the results of
Exercise \ref{exe.tridiag.cf} and Exercise \ref{exe.tridiag.fib} (the details
are left to the reader), we obtain the equality%
\begin{align*}
\dfrac{f_{n+1}}{f_{n}}  &  =1-\dfrac{1\left(  -1\right)  }{1-\dfrac{1\left(
-1\right)  }{1-\dfrac{1\left(  -1\right)  }{%
\begin{array}
[c]{ccc}%
1- &  & \\
& \ddots & \\
&  & -\dfrac{1\left(  -1\right)  }{1-\dfrac{1\left(  -1\right)  }{1}}%
\end{array}
}}}\ \ \ \ \ \ \ \ \ \ \left(  \text{with }n-1\text{ fractions in
total}\right) \\
&  =1+\dfrac{1}{1+\dfrac{1}{1+\dfrac{1}{%
\begin{array}
[c]{ccc}%
1+ &  & \\
& \ddots & \\
&  & +\dfrac{1}{1+\dfrac{1}{1}}%
\end{array}
}}}\ \ \ \ \ \ \ \ \ \ \left(  \text{with }n-1\text{ fractions in
total}\right)  .
\end{align*}
If you know
\href{https://en.wikipedia.org/wiki/Golden_ratio#Alternative_forms}{some
trivia about the golden ratio}, you might recognize this as a part of the
continued fraction for the golden ratio $\varphi$. The whole continued
fraction for $\varphi$ is%
\[
\varphi=1+\dfrac{1}{1+\dfrac{1}{1+\dfrac{1}{%
\begin{array}
[c]{cc}%
1+ & \\
& \ddots
\end{array}
}}}\ \ \ \ \ \ \ \ \ \ \left(  \text{with infinitely many fractions}\right)
.
\]
This hints at the fact that $\lim\limits_{n\rightarrow\infty}\dfrac{f_{n+1}%
}{f_{n}}=\varphi$. (This is easy to prove without continued fractions, of course.)
\end{remark}

\subsection{On block-triangular matrices}

\begin{definition}
\label{def.block2x2}Let $n$, $n^{\prime}$, $m$ and $m^{\prime}$ be four
nonnegative integers.

Let $A=\left(  a_{i,j}\right)  _{1\leq i\leq n,\ 1\leq j\leq m}$ be an
$n\times m$-matrix.

Let $B=\left(  b_{i,j}\right)  _{1\leq i\leq n,\ 1\leq j\leq m^{\prime}}$ be
an $n\times m^{\prime}$-matrix.

Let $C=\left(  c_{i,j}\right)  _{1\leq i\leq n^{\prime},\ 1\leq j\leq m}$ be
an $n^{\prime}\times m$-matrix.

Let $D=\left(  d_{i,j}\right)  _{1\leq i\leq n^{\prime},\ 1\leq j\leq
m^{\prime}}$ be an $n^{\prime}\times m^{\prime}$-matrix.

Then, $\left(
\begin{array}
[c]{cc}%
A & B\\
C & D
\end{array}
\right)  $ will mean the $\left(  n+n^{\prime}\right)  \times\left(
m+m^{\prime}\right)  $-matrix
\[
\left(
\begin{array}
[c]{cccccccc}%
a_{1,1} & a_{1,2} & \cdots & a_{1,m} & b_{1,1} & b_{1,2} & \cdots &
b_{1,m^{\prime}}\\
a_{2,1} & a_{2,2} & \cdots & a_{2,m} & b_{2,1} & b_{2,2} & \cdots &
b_{2,m^{\prime}}\\
\vdots & \vdots & \ddots & \vdots & \vdots & \vdots & \ddots & \vdots\\
a_{n,1} & a_{n,2} & \cdots & a_{n,m} & b_{n,1} & b_{n,2} & \cdots &
b_{n,m^{\prime}}\\
c_{1,1} & c_{1,2} & \cdots & c_{1,m} & d_{1,1} & d_{1,2} & \cdots &
d_{1,m^{\prime}}\\
c_{2,1} & c_{2,2} & \cdots & c_{2,m} & d_{2,1} & d_{2,2} & \cdots &
d_{2,m^{\prime}}\\
\vdots & \vdots & \ddots & \vdots & \vdots & \vdots & \ddots & \vdots\\
c_{n^{\prime},1} & c_{n^{\prime},2} & \cdots & c_{n^{\prime},m} &
d_{n^{\prime},1} & d_{n^{\prime},2} & \cdots & d_{n^{\prime},m^{\prime}}%
\end{array}
\right)  .
\]
(Formally speaking, this means that%
\begin{equation}
\left(
\begin{array}
[c]{cc}%
A & B\\
C & D
\end{array}
\right)  =\left(
\begin{cases}
a_{i,j}, & \text{if }i\leq n\text{ and }j\leq m;\\
b_{i,j-m}, & \text{if }i\leq n\text{ and }j>m;\\
c_{i-n,j}, & \text{if }i>n\text{ and }j\leq m;\\
d_{i-n,j-m}, & \text{if }i>n\text{ and }j>m
\end{cases}
\right)  _{1\leq i\leq n+n^{\prime},\ 1\leq j\leq m+m^{\prime}}.
\label{eq.def.block2x2.formal}%
\end{equation}
Less formally, we can say that $\left(
\begin{array}
[c]{cc}%
A & B\\
C & D
\end{array}
\right)  $ is the matrix obtained by gluing the matrices $A$, $B$, $C$ and $D$
to form one big $\left(  n+n^{\prime}\right)  \times\left(  m+m^{\prime
}\right)  $-matrix, where the right border of $A$ is glued together with the
left border of $B$, the bottom border of $A$ is glued together with the top
border of $C$, etc.)

Do not get fooled by the notation $\left(
\begin{array}
[c]{cc}%
A & B\\
C & D
\end{array}
\right)  $: It is (in general) not a $2\times2$-matrix, but an $\left(
n+n^{\prime}\right)  \times\left(  m+m^{\prime}\right)  $-matrix, and its
entries are not $A$, $B$, $C$ and $D$ but the entries of $A$, $B$, $C$ and $D$.
\end{definition}

\begin{example}
If $A=\left(
\begin{array}
[c]{cc}%
a_{1,1} & a_{1,2}\\
a_{2,1} & a_{2,2}%
\end{array}
\right)  $, $B=\left(
\begin{array}
[c]{c}%
b_{1}\\
b_{2}%
\end{array}
\right)  $, $C=\left(
\begin{array}
[c]{cc}%
c_{1} & c_{2}%
\end{array}
\right)  $ and $D=\left(
\begin{array}
[c]{c}%
d
\end{array}
\right)  $, then $\left(
\begin{array}
[c]{cc}%
A & B\\
C & D
\end{array}
\right)  =\left(
\begin{array}
[c]{ccc}%
a_{1,1} & a_{1,2} & b_{1}\\
a_{2,1} & a_{2,2} & b_{2}\\
c_{1} & c_{2} & d
\end{array}
\right)  $.
\end{example}

The notation $\left(
\begin{array}
[c]{cc}%
A & B\\
C & D
\end{array}
\right)  $ introduced in Definition \ref{def.block2x2} is a particular case of
a more general notation -- the \textit{block-matrix construction} -- for
gluing together multiple matrices with matching dimensions\footnote{This
construction defines an $\left(  n_{1}+n_{2}+\cdots+n_{x}\right)
\times\left(  m_{1}+m_{2}+\cdots+m_{y}\right)  $-matrix%
\begin{equation}
\left(
\begin{array}
[c]{cccc}%
A_{1,1} & A_{1,2} & \cdots & A_{1,y}\\
A_{2,1} & A_{2,2} & \cdots & A_{2,y}\\
\vdots & \vdots & \ddots & \vdots\\
A_{x,1} & A_{x,2} & \cdots & A_{x,y}%
\end{array}
\right)  \label{eq.block-general}%
\end{equation}
whenever you have given two nonnegative integers $x$ and $y$, an $x$-tuple
$\left(  n_{1},n_{2},\ldots,n_{x}\right)  \in\mathbb{N}^{x}$, a $y$-tuple
$\left(  m_{1},m_{2},\ldots,m_{y}\right)  \in\mathbb{N}^{y}$, and an
$n_{i}\times m_{j}$-matrix $A_{i,j}$ for every $i\in\left\{  1,2,\ldots
,x\right\}  $ and every $j\in\left\{  1,2,\ldots,y\right\}  $. I guess you can
guess the definition of this matrix. So you start with an \textquotedblleft%
$x\times y$-matrix of matrices\textquotedblright\ and glue them together to an
$\left(  n_{1}+n_{2}+\cdots+n_{x}\right)  \times\left(  m_{1}+m_{2}%
+\cdots+m_{y}\right)  $-matrix (provided that the dimensions of these matrices
allow them to be glued -- e.g., you cannot glue a $2\times3$-matrix to a
$4\times6$-matrix along its right border, nor on any other border).
\par
It is called \textquotedblleft block-matrix construction\textquotedblright%
\ because the original matrices $A_{i,j}$ appear as \textquotedblleft
blocks\textquotedblright\ in the big matrix (\ref{eq.block-general}). Most
authors define block matrices to be matrices which are \textquotedblleft
partitioned\textquotedblright\ into blocks as in (\ref{eq.block-general});
this is essentially our construction in reverse: Instead of gluing several
\textquotedblleft small\textquotedblright\ matrices into a big one, they study
big matrices partitioned into many small matrices. Of course, the properties
of their \textquotedblleft block matrices\textquotedblright\ are equivalent to
those of our \textquotedblleft block-matrix construction\textquotedblright.}.
We shall only need the particular case that is Definition \ref{def.block2x2}, however.

\begin{definition}
Let $n\in\mathbb{N}$ and $m\in\mathbb{N}$. Recall that $\mathbb{K}^{n\times
m}$ is the set of all $n\times m$-matrices.

We use $0_{n\times m}$ (or sometimes just $0$) to denote the $n\times
m$\textit{ zero matrix}. (As we recall, this is the $n\times m$-matrix whose
all entries are $0$; in other words, this is the $n\times m$-matrix $\left(
0\right)  _{1\leq i\leq n,\ 1\leq j\leq m}$.)
\end{definition}

\begin{exercise}
\label{exe.block2x2.mult}Let $n$, $n^{\prime}$, $m$, $m^{\prime}$, $\ell$ and
$\ell^{\prime}$ be six nonnegative integers. Let $A\in\mathbb{K}^{n\times m}$,
$B\in\mathbb{K}^{n\times m^{\prime}}$, $C\in\mathbb{K}^{n^{\prime}\times m}$,
$D\in\mathbb{K}^{n^{\prime}\times m^{\prime}}$, $A^{\prime}\in\mathbb{K}%
^{m\times\ell}$, $B^{\prime}\in\mathbb{K}^{m\times\ell^{\prime}}$, $C^{\prime
}\in\mathbb{K}^{m^{\prime}\times\ell}$ and $D^{\prime}\in\mathbb{K}%
^{m^{\prime}\times\ell^{\prime}}$. Then, prove that%
\[
\left(
\begin{array}
[c]{cc}%
A & B\\
C & D
\end{array}
\right)  \left(
\begin{array}
[c]{cc}%
A^{\prime} & B^{\prime}\\
C^{\prime} & D^{\prime}%
\end{array}
\right)  =\left(
\begin{array}
[c]{cc}%
AA^{\prime}+BC^{\prime} & AB^{\prime}+BD^{\prime}\\
CA^{\prime}+DC^{\prime} & CB^{\prime}+DD^{\prime}%
\end{array}
\right)  .
\]

\end{exercise}

\begin{remark}
The intuitive meaning of Exercise \ref{exe.block2x2.mult} is that the product
of two matrices in \textquotedblleft block-matrix notation\textquotedblright%
\ can be computed by applying the usual multiplication rule \textquotedblleft
on the level of blocks\textquotedblright, without having to fall back to
multiplying single entries. However, when applying Exercise
\ref{exe.block2x2.mult}, do not forget to check that its conditions are
satisfied. Let me give an example and a non-example:

\textbf{Example:} If $A=\left(
\begin{array}
[c]{c}%
a_{1}\\
a_{2}\\
a_{3}%
\end{array}
\right)  $, $B=\left(
\begin{array}
[c]{cc}%
b_{1,1} & b_{1,2}\\
b_{2,1} & b_{2,2}\\
b_{3,1} & b_{3,2}%
\end{array}
\right)  $, $C=\left(
\begin{array}
[c]{c}%
c
\end{array}
\right)  $, $D=\left(
\begin{array}
[c]{cc}%
d_{1} & d_{2}%
\end{array}
\right)  $, $A^{\prime}=\left(
\begin{array}
[c]{cc}%
a_{1}^{\prime} & a_{2}^{\prime}%
\end{array}
\right)  $, $B^{\prime}=\left(
\begin{array}
[c]{cc}%
b_{1}^{\prime} & b_{2}^{\prime}%
\end{array}
\right)  $, $C^{\prime}=\left(
\begin{array}
[c]{cc}%
c_{1,1}^{\prime} & c_{1,2}^{\prime}\\
c_{2,1}^{\prime} & c_{2,2}^{\prime}%
\end{array}
\right)  $ and $D^{\prime}=\left(
\begin{array}
[c]{cc}%
d_{1,1}^{\prime} & d_{1,2}^{\prime}\\
d_{2,1}^{\prime} & d_{2,2}^{\prime}%
\end{array}
\right)  $, then Exercise \ref{exe.block2x2.mult} can be applied (with $n=3$,
$n^{\prime}=1$, $m=1$, $m^{\prime}=2$, $\ell=2$ and $\ell^{\prime}=2$), and
thus we obtain%
\[
\left(
\begin{array}
[c]{cc}%
A & B\\
C & D
\end{array}
\right)  \left(
\begin{array}
[c]{cc}%
A^{\prime} & B^{\prime}\\
C^{\prime} & D^{\prime}%
\end{array}
\right)  =\left(
\begin{array}
[c]{cc}%
AA^{\prime}+BC^{\prime} & AB^{\prime}+BD^{\prime}\\
CA^{\prime}+DC^{\prime} & CB^{\prime}+DD^{\prime}%
\end{array}
\right)  .
\]

\textbf{Non-example:} If $A=\left(
\begin{array}
[c]{c}%
a_{1}\\
a_{2}%
\end{array}
\right)  $, $B=\left(
\begin{array}
[c]{cc}%
b_{1,1} & b_{1,2}\\
b_{2,1} & b_{2,2}%
\end{array}
\right)  $, $C=\left(
\begin{array}
[c]{c}%
c_{1}\\
c_{2}%
\end{array}
\right)  $, $D=\left(
\begin{array}
[c]{cc}%
d_{1,1} & d_{1,2}\\
d_{2,1} & d_{2,2}%
\end{array}
\right)  $, $A^{\prime}=\left(
\begin{array}
[c]{cc}%
a_{1,1}^{\prime} & a_{1,2}^{\prime}\\
a_{2,1}^{\prime} & a_{2,2}^{\prime}%
\end{array}
\right)  $, $B^{\prime}=\left(
\begin{array}
[c]{cc}%
b_{1,1}^{\prime} & b_{1,2}^{\prime}\\
b_{2,1}^{\prime} & b_{2,2}^{\prime}%
\end{array}
\right)  $, $C^{\prime}=\left(
\begin{array}
[c]{cc}%
c_{1}^{\prime} & c_{2}^{\prime}%
\end{array}
\right)  $ and $D^{\prime}=\left(
\begin{array}
[c]{cc}%
d_{1}^{\prime} & d_{2}^{\prime}%
\end{array}
\right)  $, then Exercise \ref{exe.block2x2.mult} cannot be applied, because
there exist no $n,m,\ell\in\mathbb{N}$ such that $A\in\mathbb{K}^{n\times m}$
and $A^{\prime}\in\mathbb{K}^{m\times\ell}$. (Indeed, the number of columns of
$A$ does not equal the number of rows of $A^{\prime}$, but these numbers would
both have to be $m$.) The matrices $\left(
\begin{array}
[c]{cc}%
A & B\\
C & D
\end{array}
\right)  $ and $\left(
\begin{array}
[c]{cc}%
A^{\prime} & B^{\prime}\\
C^{\prime} & D^{\prime}%
\end{array}
\right)  $ still exist in this case, and can even be multiplied, but their
product is not given by a simple formula such as the one in Exercise
\ref{exe.block2x2.mult}. Thus, beware of seeing Exercise
\ref{exe.block2x2.mult} as a panacea for multiplying matrices blockwise.
\end{remark}

\begin{exercise}
\label{exe.block2x2.tridet}Let $n\in\mathbb{N}$ and $m\in\mathbb{N}$. Let $A$
be an $n\times n$-matrix. Let $B$ be an $n\times m$-matrix. Let $D$ be an
$m\times m$-matrix. Prove that%
\[
\det\left(
\begin{array}
[c]{cc}%
A & B\\
0_{m\times n} & D
\end{array}
\right)  =\det A\cdot\det D.
\]

\end{exercise}

\begin{example}
Exercise \ref{exe.block2x2.tridet} (applied to $n=2$ and $m=3$) yields%
\[
\det\left(
\begin{array}
[c]{ccccc}%
a_{1,1} & a_{1,2} & b_{1,1} & b_{1,2} & b_{1,3}\\
a_{2,1} & a_{2,2} & b_{2,1} & b_{2,2} & b_{2,3}\\
0 & 0 & c_{1,1} & c_{1,2} & c_{1,3}\\
0 & 0 & c_{2,1} & c_{2,2} & c_{2,3}\\
0 & 0 & c_{3,1} & c_{3,2} & c_{3,3}%
\end{array}
\right)  =\det\left(
\begin{array}
[c]{cc}%
a_{1,1} & a_{1,2}\\
a_{2,1} & a_{2,2}%
\end{array}
\right)  \cdot\det\left(
\begin{array}
[c]{ccc}%
c_{1,1} & c_{1,2} & c_{1,3}\\
c_{2,1} & c_{2,2} & c_{2,3}\\
c_{3,1} & c_{3,2} & c_{3,3}%
\end{array}
\right)  .
\]

\end{example}

\begin{remark}
Not every determinant of the form $\det\left(
\begin{array}
[c]{cc}%
A & B\\
0_{m\times n} & D
\end{array}
\right)  $ can be computed using Exercise \ref{exe.block2x2.tridet}. In fact,
Exercise \ref{exe.block2x2.tridet} requires $A$ to be an $n\times n$-matrix
and $D$ to be an $m\times m$-matrix; thus, both $A$ and $D$ have to be square
matrices in order for Exercise \ref{exe.block2x2.tridet} to be applicable. For
instance, Exercise \ref{exe.block2x2.tridet} cannot be applied to compute
$\det\left(
\begin{array}
[c]{ccc}%
a_{1} & b_{1,1} & b_{1,2}\\
a_{2} & b_{2,1} & b_{2,2}\\
0 & c_{1} & c_{2}%
\end{array}
\right)  $.
\end{remark}

\begin{remark}
\label{rmk.block2x2.det-gen}You might wonder whether Exercise
\ref{exe.block2x2.tridet} generalizes to a formula for $\det\left(
\begin{array}
[c]{cc}%
A & B\\
C & D
\end{array}
\right)  $ when $A\in\mathbb{K}^{n\times n}$, $B\in\mathbb{K}^{n\times m}$,
$C\in\mathbb{K}^{m\times n}$ and $D\in\mathbb{K}^{m\times m}$. The general
answer is \textquotedblleft No\textquotedblright. However, when $D$ is
invertible, there exists such a formula (the
\href{https://en.wikipedia.org/wiki/Schur_complement}{Schur complement}
formula shown in Exercise \ref{exe.block2x2.schur} below). Curiously, there is
also a formula for the case when $n=m$ and $CD=DC$ (see \cite[Theorem
3]{Silvest}).
\end{remark}

We notice that Exercise \ref{exe.block2x2.tridet} allows us to solve Exercise
\ref{exe.ps4.5} in a new way.

An analogue of Exercise \ref{exe.block2x2.tridet} exists in which the
$0_{m\times n}$ in the lower-left part of the matrix is replaced by a
$0_{n\times m}$ in the upper-right part:

\begin{exercise}
\label{exe.block2x2.tridet.transposed}Let $n\in\mathbb{N}$ and $m\in
\mathbb{N}$. Let $A$ be an $n\times n$-matrix. Let $C$ be an $m\times
n$-matrix. Let $D$ be an $m\times m$-matrix. Prove that%
\[
\det\left(
\begin{array}
[c]{cc}%
A & 0_{n\times m}\\
C & D
\end{array}
\right)  =\det A\cdot\det D.
\]

\end{exercise}

\begin{exercise}
\label{exe.det.anotherpattern}\textbf{(a)} Compute the determinant of the
$7\times7$-matrix%
\[
\left(
\begin{array}
[c]{ccccccc}%
a & 0 & 0 & 0 & 0 & 0 & b\\
0 & a^{\prime} & 0 & 0 & 0 & b^{\prime} & 0\\
0 & 0 & a^{\prime\prime} & 0 & b^{\prime\prime} & 0 & 0\\
0 & 0 & 0 & e & 0 & 0 & 0\\
0 & 0 & c^{\prime\prime} & 0 & d^{\prime\prime} & 0 & 0\\
0 & c^{\prime} & 0 & 0 & 0 & d^{\prime} & 0\\
c & 0 & 0 & 0 & 0 & 0 & d
\end{array}
\right)  ,
\]
where $a,a^{\prime},a^{\prime\prime},b,b^{\prime},b^{\prime\prime}%
,c,c^{\prime},c^{\prime\prime},d,d^{\prime},d^{\prime\prime},e$ are elements
of $\mathbb{K}$.

\textbf{(b)} Compute the determinant of the $6\times6$-matrix%
\[
\left(
\begin{array}
[c]{cccccc}%
a & 0 & 0 & \ell & 0 & 0\\
0 & b & 0 & 0 & m & 0\\
0 & 0 & c & 0 & 0 & n\\
g & 0 & 0 & d & 0 & 0\\
0 & h & 0 & 0 & e & 0\\
0 & 0 & k & 0 & 0 & f
\end{array}
\right)  ,
\]
where $a,b,c,d,e,f,g,h,k,\ell,m,n$ are elements of $\mathbb{K}$.
\end{exercise}

\begin{exercise}
\label{exe.det.creative}Invent and solve an exercise on computing determinants.
\end{exercise}

\subsection{The adjugate matrix}

We start this section with a variation on Theorem \ref{thm.laplace.gen}:

\begin{proposition}
\label{prop.laplace.0}Let $n\in\mathbb{N}$. Let $A=\left(  a_{i,j}\right)
_{1\leq i\leq n,\ 1\leq j\leq n}$ be an $n\times n$-matrix. Let $r\in\left\{
1,2,\ldots,n\right\}  $.

\textbf{(a)} For every $p\in\left\{  1,2,\ldots,n\right\}  $ satisfying $p\neq
r$, we have%
\[
0=\sum_{q=1}^{n}\left(  -1\right)  ^{p+q}a_{r,q}\det\left(  A_{\sim p,\sim
q}\right)  .
\]

\textbf{(b)} For every $q\in\left\{  1,2,\ldots,n\right\}  $ satisfying $q\neq
r$, we have%
\[
0=\sum_{p=1}^{n}\left(  -1\right)  ^{p+q}a_{p,r}\det\left(  A_{\sim p,\sim
q}\right)  .
\]

\end{proposition}

\begin{vershort}
\begin{proof}
[Proof of Proposition \ref{prop.laplace.0}.]\textbf{(a)} Let $p\in\left\{
1,2,\ldots,n\right\}  $ be such that $p\neq r$.

Let $C$ be the $n\times n$-matrix obtained from $A$ by replacing the $p$-th
row of $A$ by the $r$-th row of $A$. Thus, the $p$-th and the $r$-th rows of
$C$ are equal. Therefore, the matrix $C$ has two equal rows (since $p\neq r$).
Hence, $\det C=0$ (by Exercise \ref{exe.ps4.6} \textbf{(e)}, applied to $C$
instead of $A$).

Let us write the $n\times n$-matrix $C$ in the form $C=\left(  c_{i,j}\right)
_{1\leq i\leq n,\ 1\leq j\leq n}$.

The $p$-th row of $C$ equals the $r$-th row of $A$ (by the construction of
$C$). In other words,%
\begin{equation}
c_{p,q}=a_{r,q}\ \ \ \ \ \ \ \ \ \ \text{for every }q\in\left\{
1,2,\ldots,n\right\}  . \label{pf.prop.laplace.0.short.cpq}%
\end{equation}
On the other hand, the matrix $C$ equals the matrix $A$ in all rows but the
$p$-th one (again, by the construction of $C$). Hence, if we cross out the
$p$-th rows in both $C$ and $A$, then the matrices $C$ and $A$ become equal.
Therefore,%
\begin{equation}
C_{\sim p,\sim q}=A_{\sim p,\sim q}\ \ \ \ \ \ \ \ \ \ \text{for every }%
q\in\left\{  1,2,\ldots,n\right\}  \label{pf.prop.laplace.0.short.CA}%
\end{equation}
(because the construction of $C_{\sim p,\sim q}$ from $C$ involves crossing
out the $p$-th row, and so does the construction of $A_{\sim p,\sim q}$ from
$A$).

Now, $\det C=0$, so that%
\begin{align*}
0  &  =\det C=\sum_{q=1}^{n}\left(  -1\right)  ^{p+q}\underbrace{c_{p,q}%
}_{\substack{=a_{r,q}\\\text{(by (\ref{pf.prop.laplace.0.short.cpq}))}}%
}\det\left(  \underbrace{C_{\sim p,\sim q}}_{\substack{=A_{\sim p,\sim
q}\\\text{(by (\ref{pf.prop.laplace.0.short.CA}))}}}\right) \\
&  \ \ \ \ \ \ \ \ \ \ \left(  \text{by Theorem \ref{thm.laplace.gen}
\textbf{(a)}, applied to }C\text{ and }c_{i,j}\text{ instead of }A\text{ and
}a_{i,j}\right) \\
&  =\sum_{q=1}^{n}\left(  -1\right)  ^{p+q}a_{r,q}\det\left(  A_{\sim p,\sim
q}\right)  .
\end{align*}
This proves Proposition \ref{prop.laplace.0} \textbf{(a)}.

\textbf{(b)} This proof is rather similar to the proof of Proposition
\ref{prop.laplace.0} \textbf{(a)}, except that rows are now replaced by
columns. We leave the details to the reader.
\end{proof}
\end{vershort}

\begin{verlong}
\begin{proof}
[Proof of Proposition \ref{prop.laplace.0}.]\textbf{(a)} Let $p\in\left\{
1,2,\ldots,n\right\}  $ be such that $p\neq r$.

Let $w$ be the $r$-th row of $A$ (regarded, as usual, as a row vector). Thus,
$w=\left(  \text{the }r\text{-th row of }A\right)  $.

Let $C$ be the $n\times n$-matrix obtained from $A$ by replacing the $p$-th
row of $A$ by the row vector $w$. Thus,%
\begin{align}
&  \left(  \left(  \text{the }u\text{-th row of }C\right)  =\left(  \text{the
}u\text{-th row of }A\right)  \right. \label{pf.prop.laplace.0.Cu}\\
&  \ \ \ \ \ \ \ \ \ \ \left.  \text{for all }u\in\left\{  1,2,\ldots
,n\right\}  \text{ satisfying }u\neq p\right)  ,\nonumber
\end{align}
whereas%
\begin{equation}
\left(  \text{the }p\text{-th row of }C\right)  =w.
\label{pf.prop.laplace.0.Cp}%
\end{equation}
The matrix $C$ has two equal rows\footnote{\textit{Proof.} We have $r\neq p$
(since $p\neq r$). Hence, (\ref{pf.prop.laplace.0.Cu}) (applied to $u=r$)
yields%
\begin{align*}
\left(  \text{the }r\text{-th row of }C\right)   &  =\left(  \text{the
}r\text{-th row of }A\right)  =w\ \ \ \ \ \ \ \ \ \ \left(  \text{since
}w=\left(  \text{the }r\text{-th row of }A\right)  \right) \\
&  =\left(  \text{the }p\text{-th row of }C\right)
\ \ \ \ \ \ \ \ \ \ \left(  \text{by (\ref{pf.prop.laplace.0.Cp})}\right)  .
\end{align*}
In other words, the $r$-th row of $C$ and the $p$-th row of $C$ are equal.
Since $r\neq p$, this shows that the matrix $C$ has two equal rows. Qed.}.
Hence, $\det C=0$ (by Exercise \ref{exe.ps4.6} \textbf{(e)}, applied to $C$
instead of $A$).

Let us write the $n\times n$-matrix $C$ in the form $C=\left(  c_{i,j}\right)
_{1\leq i\leq n,\ 1\leq j\leq n}$. Thus, for every $u\in\left\{
1,2,\ldots,n\right\}  $, we have%
\begin{equation}
\left(  \text{the }u\text{-th row of }C\right)  =\left(  c_{u,1}%
,c_{u,2},\ldots,c_{u,n}\right)  . \label{pf.prop.laplace.0.u-th-row}%
\end{equation}
Applying this to $u=p$, we obtain%
\[
\left(  \text{the }p\text{-th row of }C\right)  =\left(  c_{p,1}%
,c_{p,2},\ldots,c_{p,n}\right)  ,
\]
so that%
\begin{align*}
\left(  c_{p,1},c_{p,2},\ldots,c_{p,n}\right)   &  =\left(  \text{the
}p\text{-th row of }C\right)  =w=\left(  \text{the }r\text{-th row of
}A\right) \\
&  =\left(  a_{r,1},a_{r,2},\ldots,a_{r,n}\right)  \ \ \ \ \ \ \ \ \ \ \left(
\text{since }A=\left(  a_{i,j}\right)  _{1\leq i\leq n,\ 1\leq j\leq
n}\right)  .
\end{align*}
In other words,%
\begin{equation}
c_{p,q}=a_{r,q}\ \ \ \ \ \ \ \ \ \ \text{for every }q\in\left\{
1,2,\ldots,n\right\}  . \label{pf.prop.laplace.0.cpq}%
\end{equation}

On the other hand,%
\begin{equation}
c_{u,q}=a_{u,q}\ \ \ \ \ \ \ \ \ \ \text{for every }q\in\left\{
1,2,\ldots,n\right\}  \text{ and }u\in\left\{  1,2,\ldots,n\right\}  \text{
satisfying }u\neq p \label{pf.prop.laplace.0.cuq}%
\end{equation}
\footnote{\textit{Proof of (\ref{pf.prop.laplace.0.cuq}):} Let $u\in\left\{
1,2,\ldots,n\right\}  $ be such that $u\neq p$. Thus,%
\begin{align*}
\left(  c_{u,1},c_{u,2},\ldots,c_{u,n}\right)   &  =\left(  \text{the
}u\text{-th row of }C\right)  \ \ \ \ \ \ \ \ \ \ \left(  \text{by
(\ref{pf.prop.laplace.0.u-th-row})}\right) \\
&  =\left(  \text{the }u\text{-th row of }A\right)
\ \ \ \ \ \ \ \ \ \ \left(  \text{by (\ref{pf.prop.laplace.0.Cu})}\right) \\
&  =\left(  a_{u,1},a_{u,2},\ldots,a_{u,n}\right)  \ \ \ \ \ \ \ \ \ \ \left(
\text{since }A=\left(  a_{i,j}\right)  _{1\leq i\leq n,\ 1\leq j\leq
n}\right)  .
\end{align*}
In other words, $c_{u,q}=a_{u,q}$ for every $q\in\left\{  1,2,\ldots
,n\right\}  $. This proves (\ref{pf.prop.laplace.0.cuq}).}. Now, it is easy to
see that%
\begin{equation}
C_{\sim p,\sim q}=A_{\sim p,\sim q}\ \ \ \ \ \ \ \ \ \ \text{for every }%
q\in\left\{  1,2,\ldots,n\right\}  \label{pf.prop.laplace.0.CA}%
\end{equation}
\footnote{\textit{Proof of (\ref{pf.prop.laplace.0.CA}):} Let $q\in\left\{
1,2,\ldots,n\right\}  $.
\par
Let $\left(  u_{1},u_{2},\ldots,u_{n-1}\right)  $ denote the $\left(
n-1\right)  $-tuple $\left(  1,2,\ldots,\widehat{p},\ldots,n\right)  $. Thus,
$\left(  u_{1},u_{2},\ldots,u_{n-1}\right)  =\left(  1,2,\ldots,\widehat{p}%
,\ldots,n\right)  $, so that $\left\{  u_{1},u_{2},\ldots,u_{n-1}\right\}
=\left\{  1,2,\ldots,\widehat{p},\ldots,n\right\}  =\left\{  1,2,\ldots
,n\right\}  \setminus\left\{  p\right\}  $.
\par
Now, let $x\in\left\{  1,2,\ldots,n-1\right\}  $. Then, $u_{x}\in\left\{
u_{1},u_{2},\ldots,u_{n-1}\right\}  =\left\{  1,2,\ldots,n\right\}
\setminus\left\{  p\right\}  $, so that $u_{x}\neq p$. Hence,%
\begin{equation}
c_{u_{x},q}=a_{u_{x},q}\ \ \ \ \ \ \ \ \ \ \text{for every }q\in\left\{
1,2,\ldots,n\right\}  \label{pf.prop.laplace.0.CA.pf.1}%
\end{equation}
(by (\ref{pf.prop.laplace.0.cuq}), applied to $u=u_{x}$).
\par
Let us now forget that we fixed $x$. We thus have shown that
(\ref{pf.prop.laplace.0.CA.pf.1}) holds for every $x\in\left\{  1,2,\ldots
,n-1\right\}  $.
\par
Let $\left(  v_{1},v_{2},\ldots,v_{n-1}\right)  $ denote the $\left(
n-1\right)  $-tuple $\left(  1,2,\ldots,\widehat{q},\ldots,n\right)  $. Thus,
$\left(  v_{1},v_{2},\ldots,v_{n-1}\right)  =\left(  1,2,\ldots,\widehat{q}%
,\ldots,n\right)  $.
\par
Now, the definition of $C_{\sim p,\sim q}$ yields%
\begin{align*}
C_{\sim p,\sim q}  &  =\operatorname*{sub}\nolimits_{1,2,\ldots,\widehat{p}%
,\ldots,n}^{1,2,\ldots,\widehat{q},\ldots,n}C\\
&  =\operatorname*{sub}\nolimits_{1,2,\ldots,\widehat{p},\ldots,n}%
^{v_{1},v_{2},\ldots,v_{n-1}}C\ \ \ \ \ \ \ \ \ \ \left(  \text{since }\left(
1,2,\ldots,\widehat{q},\ldots,n\right)  =\left(  v_{1},v_{2},\ldots
,v_{n-1}\right)  \right) \\
&  =\operatorname*{sub}\nolimits_{u_{1},u_{2},\ldots,u_{n-1}}^{v_{1}%
,v_{2},\ldots,v_{n-1}}C\ \ \ \ \ \ \ \ \ \ \left(  \text{since }\left(
1,2,\ldots,\widehat{p},\ldots,n\right)  =\left(  u_{1},u_{2},\ldots
,u_{n-1}\right)  \right) \\
&  =\left(  \underbrace{c_{u_{x},v_{y}}}_{\substack{=a_{u_{x},v_{y}%
}\\\text{(by (\ref{pf.prop.laplace.0.CA.pf.1}),}\\\text{applied to }%
q=v_{y}\text{)}}}\right)  _{1\leq x\leq n-1,\ 1\leq y\leq n-1}\\
&  \ \ \ \ \ \ \ \ \ \ \left(  \text{by the definition of }\operatorname*{sub}%
\nolimits_{u_{1},u_{2},\ldots,u_{n-1}}^{v_{1},v_{2},\ldots,v_{n-1}}C\text{,
since }C=\left(  c_{i,j}\right)  _{1\leq i\leq n,\ 1\leq j\leq n}\right) \\
&  =\left(  a_{u_{x},v_{y}}\right)  _{1\leq x\leq n-1,\ 1\leq y\leq n-1}.
\end{align*}
Compared with%
\begin{align*}
A_{\sim p,\sim q}  &  =\operatorname*{sub}\nolimits_{1,2,\ldots,\widehat{p}%
,\ldots,n}^{1,2,\ldots,\widehat{q},\ldots,n}A\\
&  =\operatorname*{sub}\nolimits_{1,2,\ldots,\widehat{p},\ldots,n}%
^{v_{1},v_{2},\ldots,v_{n-1}}A\ \ \ \ \ \ \ \ \ \ \left(  \text{since }\left(
1,2,\ldots,\widehat{q},\ldots,n\right)  =\left(  v_{1},v_{2},\ldots
,v_{n-1}\right)  \right) \\
&  =\operatorname*{sub}\nolimits_{u_{1},u_{2},\ldots,u_{n-1}}^{v_{1}%
,v_{2},\ldots,v_{n-1}}A\ \ \ \ \ \ \ \ \ \ \left(  \text{since }\left(
1,2,\ldots,\widehat{p},\ldots,n\right)  =\left(  u_{1},u_{2},\ldots
,u_{n-1}\right)  \right) \\
&  =\left(  a_{u_{x},v_{y}}\right)  _{1\leq x\leq n-1,\ 1\leq y\leq n-1}\\
&  \ \ \ \ \ \ \ \ \ \ \left(  \text{by the definition of }\operatorname*{sub}%
\nolimits_{u_{1},u_{2},\ldots,u_{n-1}}^{v_{1},v_{2},\ldots,v_{n-1}}A\text{,
since }A=\left(  a_{i,j}\right)  _{1\leq i\leq n,\ 1\leq j\leq n}\right)  ,
\end{align*}
this yields $C_{\sim p,\sim q}=A_{\sim p,\sim q}$. This proves
(\ref{pf.prop.laplace.0.CA}).}. Now, $\det C=0$, so that%
\begin{align*}
0  &  =\det C=\sum_{q=1}^{n}\left(  -1\right)  ^{p+q}\underbrace{c_{p,q}%
}_{\substack{=a_{r,q}\\\text{(by (\ref{pf.prop.laplace.0.cpq}))}}}\det\left(
\underbrace{C_{\sim p,\sim q}}_{\substack{=A_{\sim p,\sim q}\\\text{(by
(\ref{pf.prop.laplace.0.CA}))}}}\right) \\
&  \ \ \ \ \ \ \ \ \ \ \left(  \text{by Theorem \ref{thm.laplace.gen}
\textbf{(a)}, applied to }C\text{ and }c_{i,j}\text{ instead of }A\text{ and
}a_{i,j}\right) \\
&  =\sum_{q=1}^{n}\left(  -1\right)  ^{p+q}a_{r,q}\det\left(  A_{\sim p,\sim
q}\right)  .
\end{align*}
This proves Proposition \ref{prop.laplace.0} \textbf{(a)}.

\textbf{(b)} This proof is rather similar to the proof of Proposition
\ref{prop.laplace.0} \textbf{(a)}, except that rows are now replaced by
columns. Let me nevertheless show this proof in full detail, for the sake of completeness:

Let $q\in\left\{  1,2,\ldots,n\right\}  $ be such that $q\neq r$.

Let $w$ be the $r$-th column of $A$ (regarded, as usual, as a column vector).
Thus, $w=\left(  \text{the }r\text{-th column of }A\right)  $.

Let $C$ be the $n\times n$-matrix obtained from $A$ by replacing the $q$-th
column of $A$ by the column vector $w$. Thus,%
\begin{align}
&  \left(  \left(  \text{the }u\text{-th column of }C\right)  =\left(
\text{the }u\text{-th column of }A\right)  \right.
\label{pf.prop.laplace.0.b.Cu}\\
&  \ \ \ \ \ \ \ \ \ \ \left.  \text{for all }u\in\left\{  1,2,\ldots
,n\right\}  \text{ satisfying }u\neq q\right)  ,\nonumber
\end{align}
whereas%
\begin{equation}
\left(  \text{the }q\text{-th column of }C\right)  =w.
\label{pf.prop.laplace.0.b.Cp}%
\end{equation}
The matrix $C$ has two equal columns\footnote{\textit{Proof.} We have $r\neq
q$ (since $q\neq r$). Hence, (\ref{pf.prop.laplace.0.b.Cu}) (applied to $u=r$)
yields%
\begin{align*}
\left(  \text{the }r\text{-th column of }C\right)   &  =\left(  \text{the
}r\text{-th column of }A\right)  =w\ \ \ \ \ \ \ \ \ \ \left(  \text{since
}w=\left(  \text{the }r\text{-th column of }A\right)  \right) \\
&  =\left(  \text{the }q\text{-th column of }C\right)
\ \ \ \ \ \ \ \ \ \ \left(  \text{by (\ref{pf.prop.laplace.0.b.Cp})}\right)  .
\end{align*}
In other words, the $r$-th column of $C$ and the $q$-th column of $C$ are
equal. Since $r\neq q$, this shows that the matrix $C$ has two equal columns.
Qed.}. Hence, $\det C=0$ (by Exercise \ref{exe.ps4.6} \textbf{(f)}, applied to
$C$ instead of $A$).

Let us write the $n\times n$-matrix $C$ in the form $C=\left(  c_{i,j}\right)
_{1\leq i\leq n,\ 1\leq j\leq n}$. Thus, for every $u\in\left\{
1,2,\ldots,n\right\}  $, we have%
\begin{equation}
\left(  \text{the }u\text{-th column of }C\right)  =\left(
\begin{array}
[c]{c}%
c_{1,u}\\
c_{2,u}\\
\vdots\\
c_{n,u}%
\end{array}
\right)  . \label{pf.prop.laplace.0.b.u-th-row}%
\end{equation}
Applying this to $u=q$, we obtain%
\[
\left(  \text{the }q\text{-th column of }C\right)  =\left(
\begin{array}
[c]{c}%
c_{1,q}\\
c_{2,q}\\
\vdots\\
c_{n,q}%
\end{array}
\right)  ,
\]
so that%
\begin{align*}
\left(
\begin{array}
[c]{c}%
c_{1,q}\\
c_{2,q}\\
\vdots\\
c_{n,q}%
\end{array}
\right)   &  =\left(  \text{the }q\text{-th column of }C\right)  =w=\left(
\text{the }r\text{-th column of }A\right) \\
&  =\left(
\begin{array}
[c]{c}%
a_{1,r}\\
a_{2,r}\\
\vdots\\
a_{n,r}%
\end{array}
\right)  \ \ \ \ \ \ \ \ \ \ \left(  \text{since }A=\left(  a_{i,j}\right)
_{1\leq i\leq n,\ 1\leq j\leq n}\right)  .
\end{align*}
In other words,%
\begin{equation}
c_{p,q}=a_{p,r}\ \ \ \ \ \ \ \ \ \ \text{for every }p\in\left\{
1,2,\ldots,n\right\}  . \label{pf.prop.laplace.0.b.cpq}%
\end{equation}

On the other hand,%
\begin{equation}
c_{p,u}=a_{p,u}\ \ \ \ \ \ \ \ \ \ \text{for every }p\in\left\{
1,2,\ldots,n\right\}  \text{ and }u\in\left\{  1,2,\ldots,n\right\}  \text{
satisfying }u\neq q \label{pf.prop.laplace.0.b.cuq}%
\end{equation}
\footnote{\textit{Proof of (\ref{pf.prop.laplace.0.b.cuq}):} Let $u\in\left\{
1,2,\ldots,n\right\}  $ be such that $u\neq q$. Thus,%
\begin{align*}
\left(
\begin{array}
[c]{c}%
c_{1,u}\\
c_{2,u}\\
\vdots\\
c_{n,u}%
\end{array}
\right)   &  =\left(  \text{the }u\text{-th column of }C\right)
\ \ \ \ \ \ \ \ \ \ \left(  \text{by (\ref{pf.prop.laplace.0.b.u-th-row}%
)}\right) \\
&  =\left(  \text{the }u\text{-th column of }A\right)
\ \ \ \ \ \ \ \ \ \ \left(  \text{by (\ref{pf.prop.laplace.0.b.Cu})}\right) \\
&  =\left(
\begin{array}
[c]{c}%
a_{1,u}\\
a_{2,u}\\
\vdots\\
a_{n,u}%
\end{array}
\right)  \ \ \ \ \ \ \ \ \ \ \left(  \text{since }A=\left(  a_{i,j}\right)
_{1\leq i\leq n,\ 1\leq j\leq n}\right)  .
\end{align*}
In other words, $c_{p,u}=a_{p,u}$ for every $p\in\left\{  1,2,\ldots
,n\right\}  $. This proves (\ref{pf.prop.laplace.0.b.cuq}).}. Now, it is easy
to see that%
\begin{equation}
C_{\sim p,\sim q}=A_{\sim p,\sim q}\ \ \ \ \ \ \ \ \ \ \text{for every }%
p\in\left\{  1,2,\ldots,n\right\}  \label{pf.prop.laplace.0.b.CA}%
\end{equation}
\footnote{\textit{Proof of (\ref{pf.prop.laplace.0.b.CA}):} Let $p\in\left\{
1,2,\ldots,n\right\}  $.
\par
Let $\left(  u_{1},u_{2},\ldots,u_{n-1}\right)  $ denote the $\left(
n-1\right)  $-tuple $\left(  1,2,\ldots,\widehat{q},\ldots,n\right)  $. Thus,
$\left(  u_{1},u_{2},\ldots,u_{n-1}\right)  =\left(  1,2,\ldots,\widehat{q}%
,\ldots,n\right)  $, so that $\left\{  u_{1},u_{2},\ldots,u_{n-1}\right\}
=\left\{  1,2,\ldots,\widehat{q},\ldots,n\right\}  =\left\{  1,2,\ldots
,n\right\}  \setminus\left\{  q\right\}  $.
\par
Now, let $y\in\left\{  1,2,\ldots,n-1\right\}  $. Then, $u_{y}\in\left\{
u_{1},u_{2},\ldots,u_{n-1}\right\}  =\left\{  1,2,\ldots,n\right\}
\setminus\left\{  q\right\}  $, so that $u_{y}\neq q$. Hence,%
\begin{equation}
c_{p,u_{y}}=a_{p,u_{y}}\ \ \ \ \ \ \ \ \ \ \text{for every }p\in\left\{
1,2,\ldots,n\right\}  \label{pf.prop.laplace.0.b.CA.pf.1}%
\end{equation}
(by (\ref{pf.prop.laplace.0.b.cuq}), applied to $u=u_{y}$).
\par
Let us now forget that we fixed $y$. We thus have shown that
(\ref{pf.prop.laplace.0.b.CA.pf.1}) holds for every $y\in\left\{
1,2,\ldots,n-1\right\}  $.
\par
Let $\left(  v_{1},v_{2},\ldots,v_{n-1}\right)  $ denote the $\left(
n-1\right)  $-tuple $\left(  1,2,\ldots,\widehat{p},\ldots,n\right)  $. Thus,
$\left(  v_{1},v_{2},\ldots,v_{n-1}\right)  =\left(  1,2,\ldots,\widehat{p}%
,\ldots,n\right)  $.
\par
Now, the definition of $C_{\sim p,\sim q}$ yields%
\begin{align*}
C_{\sim p,\sim q}  &  =\operatorname*{sub}\nolimits_{1,2,\ldots,\widehat{p}%
,\ldots,n}^{1,2,\ldots,\widehat{q},\ldots,n}C\\
&  =\operatorname*{sub}\nolimits_{v_{1},v_{2},\ldots,v_{n-1}}^{1,2,\ldots
,\widehat{q},\ldots,n}C\ \ \ \ \ \ \ \ \ \ \left(  \text{since }\left(
1,2,\ldots,\widehat{p},\ldots,n\right)  =\left(  v_{1},v_{2},\ldots
,v_{n-1}\right)  \right) \\
&  =\operatorname*{sub}\nolimits_{v_{1},v_{2},\ldots,v_{n-1}}^{u_{1}%
,u_{2},\ldots,u_{n-1}}C\ \ \ \ \ \ \ \ \ \ \left(  \text{since }\left(
1,2,\ldots,\widehat{q},\ldots,n\right)  =\left(  u_{1},u_{2},\ldots
,u_{n-1}\right)  \right) \\
&  =\left(  \underbrace{c_{v_{x},u_{y}}}_{\substack{=a_{v_{x},u_{y}%
}\\\text{(by (\ref{pf.prop.laplace.0.b.CA.pf.1}),}\\\text{applied to }%
p=v_{x}\text{)}}}\right)  _{1\leq x\leq n-1,\ 1\leq y\leq n-1}\\
&  \ \ \ \ \ \ \ \ \ \ \left(  \text{by the definition of }\operatorname*{sub}%
\nolimits_{v_{1},v_{2},\ldots,v_{n-1}}^{u_{1},u_{2},\ldots,u_{n-1}}C\text{,
since }C=\left(  c_{i,j}\right)  _{1\leq i\leq n,\ 1\leq j\leq n}\right) \\
&  =\left(  a_{v_{x},u_{y}}\right)  _{1\leq x\leq n-1,\ 1\leq y\leq n-1}.
\end{align*}
Compared with%
\begin{align*}
A_{\sim p,\sim q}  &  =\operatorname*{sub}\nolimits_{1,2,\ldots,\widehat{p}%
,\ldots,n}^{1,2,\ldots,\widehat{q},\ldots,n}A\\
&  =\operatorname*{sub}\nolimits_{v_{1},v_{2},\ldots,v_{n-1}}^{1,2,\ldots
,\widehat{q},\ldots,n}A\ \ \ \ \ \ \ \ \ \ \left(  \text{since }\left(
1,2,\ldots,\widehat{p},\ldots,n\right)  =\left(  v_{1},v_{2},\ldots
,v_{n-1}\right)  \right) \\
&  =\operatorname*{sub}\nolimits_{v_{1},v_{2},\ldots,v_{n-1}}^{u_{1}%
,u_{2},\ldots,u_{n-1}}A\ \ \ \ \ \ \ \ \ \ \left(  \text{since }\left(
1,2,\ldots,\widehat{q},\ldots,n\right)  =\left(  u_{1},u_{2},\ldots
,u_{n-1}\right)  \right) \\
&  =\left(  a_{v_{x},u_{y}}\right)  _{1\leq x\leq n-1,\ 1\leq y\leq n-1}\\
&  \ \ \ \ \ \ \ \ \ \ \left(  \text{by the definition of }\operatorname*{sub}%
\nolimits_{v_{1},v_{2},\ldots,v_{n-1}}^{u_{1},u_{2},\ldots,u_{n-1}}A\text{,
since }A=\left(  a_{i,j}\right)  _{1\leq i\leq n,\ 1\leq j\leq n}\right)  ,
\end{align*}
this yields $C_{\sim p,\sim q}=A_{\sim p,\sim q}$. This proves
(\ref{pf.prop.laplace.0.b.CA}).}. Now, $\det C=0$, so that%
\begin{align*}
0  &  =\det C=\sum_{p=1}^{n}\left(  -1\right)  ^{p+q}\underbrace{c_{p,q}%
}_{\substack{=a_{p,r}\\\text{(by (\ref{pf.prop.laplace.0.b.cpq}))}}%
}\det\left(  \underbrace{C_{\sim p,\sim q}}_{\substack{=A_{\sim p,\sim
q}\\\text{(by (\ref{pf.prop.laplace.0.b.CA}))}}}\right) \\
&  \ \ \ \ \ \ \ \ \ \ \left(  \text{by Theorem \ref{thm.laplace.gen}
\textbf{(b)}, applied to }C\text{ and }c_{i,j}\text{ instead of }A\text{ and
}a_{i,j}\right) \\
&  =\sum_{q=1}^{n}\left(  -1\right)  ^{p+q}a_{p,r}\det\left(  A_{\sim p,\sim
q}\right)  .
\end{align*}
This proves Proposition \ref{prop.laplace.0} \textbf{(b)}.
\end{proof}
\end{verlong}

We now can define the \textquotedblleft adjugate\textquotedblright\ of a matrix:

\begin{definition}
\label{def.adj}Let $n\in\mathbb{N}$. Let $A$ be an $n\times n$-matrix. We
define a new $n\times n$-matrix $\operatorname*{adj}A$ by%
\[
\operatorname*{adj}A=\left(  \left(  -1\right)  ^{i+j}\det\left(  A_{\sim
j,\sim i}\right)  \right)  _{1\leq i\leq n,\ 1\leq j\leq n}.
\]

This matrix $\operatorname*{adj}A$ is called the \textit{adjugate} of the
matrix $A$. (Some authors call it the \textquotedblleft
adjunct\textquotedblright\ or \textquotedblleft adjoint\textquotedblright\ or
\textquotedblleft classical adjoint\textquotedblright\ of $A$ instead.
However, beware of the word \textquotedblleft adjoint\textquotedblright: It
means too many different things; in particular it has a second meaning for a matrix.)
\end{definition}

The appearance of $A_{\sim j,\sim i}$ (not $A_{\sim i,\sim j}$) in Definition
\ref{def.adj} might be surprising, but it is not a mistake. We will soon see
what it is good for.

There is also a related notion, namely that of a \textquotedblleft cofactor
matrix\textquotedblright. The \textit{cofactor matrix} of an $n\times
n$-matrix $A$ is defined to be $\left(  \left(  -1\right)  ^{i+j}\det\left(
A_{\sim i,\sim j}\right)  \right)  _{1\leq i\leq n,\ 1\leq j\leq n}$. This is,
of course, the transpose $\left(  \operatorname*{adj}A\right)  ^{T}$ of
$\operatorname*{adj}A$. The entries of this matrix are called the
\textit{cofactors} of $A$.

\begin{example}
The adjugate of the $0\times0$-matrix is the $0\times0$-matrix.

The adjugate of a $1\times1$-matrix $\left(
\begin{array}
[c]{c}%
a
\end{array}
\right)  $ is $\operatorname*{adj}\left(
\begin{array}
[c]{c}%
a
\end{array}
\right)  =\left(
\begin{array}
[c]{c}%
1
\end{array}
\right)  $. (Yes, this shows that all $1\times1$-matrices have the same adjugate.)

The adjugate of a $2\times2$-matrix $\left(
\begin{array}
[c]{cc}%
a & b\\
c & d
\end{array}
\right)  $ is
\[
\operatorname*{adj}\left(
\begin{array}
[c]{cc}%
a & b\\
c & d
\end{array}
\right)  =\left(
\begin{array}
[c]{cc}%
d & -b\\
-c & a
\end{array}
\right)  .
\]

The adjugate of a $3\times3$-matrix $\left(
\begin{array}
[c]{ccc}%
a & b & c\\
d & e & f\\
g & h & i
\end{array}
\right)  $ is
\[
\operatorname*{adj}\left(
\begin{array}
[c]{ccc}%
a & b & c\\
d & e & f\\
g & h & i
\end{array}
\right)  =\left(
\begin{array}
[c]{ccc}%
ei-fh & ch-bi & bf-ce\\
fg-di & ai-cg & cd-af\\
dh-ge & bg-ah & ae-bd
\end{array}
\right)  .
\]

\end{example}

\begin{proposition}
\label{prop.adj.transpose}Let $n\in\mathbb{N}$. Let $A$ be an $n\times
n$-matrix. Then, $\operatorname*{adj}\left(  A^{T}\right)  =\left(
\operatorname*{adj}A\right)  ^{T}$.
\end{proposition}

\begin{proof}
[Proof of Proposition \ref{prop.adj.transpose}.]Let $i\in\left\{
1,2,\ldots,n\right\}  $ and $j\in\left\{  1,2,\ldots,n\right\}  $. \newline
From $i\in\left\{  1,2,\ldots,n\right\}  $, we obtain $1\leq i\leq n$, so that
$n\geq1$ and thus $n-1\in\mathbb{N}$.

The definition of $A_{\sim i,\sim j}$ yields $A_{\sim i,\sim j}%
=\operatorname*{sub}\nolimits_{1,2,\ldots,\widehat{i},\ldots,n}^{1,2,\ldots
,\widehat{j},\ldots,n}A$. But the definition of $\left(  A^{T}\right)  _{\sim
j,\sim i}$ yields%
\begin{equation}
\left(  A^{T}\right)  _{\sim j,\sim i}=\operatorname*{sub}%
\nolimits_{1,2,\ldots,\widehat{j},\ldots,n}^{1,2,\ldots,\widehat{i},\ldots
,n}\left(  A^{T}\right)  . \label{pf.prop.adj.transpose.1}%
\end{equation}

On the other hand, Proposition \ref{prop.submatrix.easy} \textbf{(e)} (applied
to $m=n$, $u=n-1$, $v=n-1$, $\left(  i_{1},i_{2},\ldots,i_{u}\right)  =\left(
1,2,\ldots,\widehat{i},\ldots,n\right)  $ and $\left(  j_{1},j_{2}%
,\ldots,j_{v}\right)  =\left(  1,2,\ldots,\widehat{j},\ldots,n\right)  $)
yields $\left(  \operatorname*{sub}\nolimits_{1,2,\ldots,\widehat{i},\ldots
,n}^{1,2,\ldots,\widehat{j},\ldots,n}A\right)  ^{T}=\operatorname*{sub}%
\nolimits_{1,2,\ldots,\widehat{j},\ldots,n}^{1,2,\ldots,\widehat{i},\ldots
,n}\left(  A^{T}\right)  $. Compared with (\ref{pf.prop.adj.transpose.1}),
this yields%
\[
\left(  A^{T}\right)  _{\sim j,\sim i}=\left(  \underbrace{\operatorname*{sub}%
\nolimits_{1,2,\ldots,\widehat{i},\ldots,n}^{1,2,\ldots,\widehat{j},\ldots
,n}A}_{=A_{\sim i,\sim j}}\right)  ^{T}=\left(  A_{\sim i,\sim j}\right)
^{T}.
\]
Hence,%
\begin{equation}
\det\left(  \underbrace{\left(  A^{T}\right)  _{\sim j,\sim i}}_{=\left(
A_{\sim i,\sim j}\right)  ^{T}}\right)  =\det\left(  \left(  A_{\sim i,\sim
j}\right)  ^{T}\right)  =\det\left(  A_{\sim i,\sim j}\right)
\label{pf.prop.adj.transpose.4}%
\end{equation}
(by Exercise \ref{exe.ps4.4}, applied to $n-1$ and $A_{\sim i,\sim j}$ instead
of $n$ and $A$).

Let us now forget that we fixed $i$ and $j$. We thus have shown that
(\ref{pf.prop.adj.transpose.4}) holds for every $i\in\left\{  1,2,\ldots
,n\right\}  $ and $j\in\left\{  1,2,\ldots,n\right\}  $.

Now, $\operatorname*{adj}A=\left(  \left(  -1\right)  ^{i+j}\det\left(
A_{\sim j,\sim i}\right)  \right)  _{1\leq i\leq n,\ 1\leq j\leq n}$, and thus
the definition of the transpose of a matrix shows that%
\[
\left(  \operatorname*{adj}A\right)  ^{T}=\left(  \underbrace{\left(
-1\right)  ^{j+i}}_{=\left(  -1\right)  ^{i+j}}\det\left(  A_{\sim i,\sim
j}\right)  \right)  _{1\leq i\leq n,\ 1\leq j\leq n}=\left(  \left(
-1\right)  ^{i+j}\det\left(  A_{\sim i,\sim j}\right)  \right)  _{1\leq i\leq
n,\ 1\leq j\leq n}.
\]
Compared with%
\begin{align*}
\operatorname*{adj}\left(  A^{T}\right)   &  =\left(  \left(  -1\right)
^{i+j}\underbrace{\det\left(  \left(  A^{T}\right)  _{\sim j,\sim i}\right)
}_{\substack{=\det\left(  A_{\sim i,\sim j}\right)  \\\text{(by
(\ref{pf.prop.adj.transpose.4}))}}}\right)  _{1\leq i\leq n,\ 1\leq j\leq n}\\
&  \ \ \ \ \ \ \ \ \ \ \left(  \text{by the definition of }\operatorname*{adj}%
\left(  A^{T}\right)  \right) \\
&  =\left(  \left(  -1\right)  ^{i+j}\det\left(  A_{\sim i,\sim j}\right)
\right)  _{1\leq i\leq n,\ 1\leq j\leq n},
\end{align*}
this yields $\operatorname*{adj}\left(  A^{T}\right)  =\left(
\operatorname*{adj}A\right)  ^{T}$. This proves Proposition
\ref{prop.adj.transpose}.
\end{proof}

The most important property of adjugates, however, is the following fact:

\begin{theorem}
\label{thm.adj.inverse}Let $n\in\mathbb{N}$. Let $A$ be an $n\times n$-matrix.
Then,%
\[
A\cdot\operatorname*{adj}A=\operatorname*{adj}A\cdot A=\det A\cdot I_{n}.
\]
(Recall that $I_{n}$ denotes the $n\times n$ identity matrix. Expressions such
as $\operatorname*{adj}A\cdot A$ and $\det A\cdot I_{n}$ have to be understood
as $\left(  \operatorname*{adj}A\right)  \cdot A$ and $\left(  \det A\right)
\cdot I_{n}$, respectively.)
\end{theorem}

\begin{example}
Recall that the adjugate of a $2\times2$-matrix is given by the formula
$\operatorname*{adj}\left(
\begin{array}
[c]{cc}%
a & b\\
c & d
\end{array}
\right)  =\left(
\begin{array}
[c]{cc}%
d & -b\\
-c & a
\end{array}
\right)  $. Thus, Theorem \ref{thm.adj.inverse} (applied to $n=2$) yields%
\[
\left(
\begin{array}
[c]{cc}%
a & b\\
c & d
\end{array}
\right)  \cdot\left(
\begin{array}
[c]{cc}%
d & -b\\
-c & a
\end{array}
\right)  =\left(
\begin{array}
[c]{cc}%
d & -b\\
-c & a
\end{array}
\right)  \cdot\left(
\begin{array}
[c]{cc}%
a & b\\
c & d
\end{array}
\right)  =\det\left(
\begin{array}
[c]{cc}%
a & b\\
c & d
\end{array}
\right)  \cdot I_{2}.
\]
(Of course, $\det\left(
\begin{array}
[c]{cc}%
a & b\\
c & d
\end{array}
\right)  \cdot I_{2}=\left(  ad-bc\right)  \cdot I_{2}=\left(
\begin{array}
[c]{cc}%
ad-bc & 0\\
0 & ad-bc
\end{array}
\right)  $.)
\end{example}

\begin{vershort}
\begin{proof}
[Proof of Theorem \ref{thm.adj.inverse}.]For any two objects $i$ and $j$, we
define $\delta_{i,j}$ to be the element $%
\begin{cases}
1, & \text{if }i=j;\\
0, & \text{if }i\neq j
\end{cases}
$ of $\mathbb{K}$. Then, $I_{n}=\left(  \delta_{i,j}\right)  _{1\leq i\leq
n,\ 1\leq j\leq n}$ (by the definition of $I_{n}$), and thus%
\begin{equation}
\det A\cdot\underbrace{I_{n}}_{=\left(  \delta_{i,j}\right)  _{1\leq i\leq
n,\ 1\leq j\leq n}}=\det A\cdot\left(  \delta_{i,j}\right)  _{1\leq i\leq
n,\ 1\leq j\leq n}=\left(  \det A\cdot\delta_{i,j}\right)  _{1\leq i\leq
n,\ 1\leq j\leq n}. \label{pf.thm.adj.inverse.short.R}%
\end{equation}

On the other hand, let us write the matrix $A$ in the form $A=\left(
a_{i,j}\right)  _{1\leq i\leq n,\ 1\leq j\leq n}$. Then, the definition of the
product of two matrices shows that%
\begin{align}
&  A\cdot\operatorname*{adj}A\nonumber\\
&  =\left(  \sum_{k=1}^{n}a_{i,k}\left(  -1\right)  ^{k+j}\det\left(  A_{\sim
j,\sim k}\right)  \right)  _{1\leq i\leq n,\ 1\leq j\leq n}\nonumber\\
&  \ \ \ \ \ \ \ \ \ \ \left(
\begin{array}
[c]{c}%
\text{since }A=\left(  a_{i,j}\right)  _{1\leq i\leq n,\ 1\leq j\leq n}\\
\text{and }\operatorname*{adj}A=\left(  \left(  -1\right)  ^{i+j}\det\left(
A_{\sim j,\sim i}\right)  \right)  _{1\leq i\leq n,\ 1\leq j\leq n}%
\end{array}
\right) \nonumber\\
&  =\left(  \sum_{q=1}^{n}\underbrace{a_{i,q}\left(  -1\right)  ^{q+j}%
}_{=\left(  -1\right)  ^{q+j}a_{i,q}}\det\left(  A_{\sim j,\sim q}\right)
\right)  _{1\leq i\leq n,\ 1\leq j\leq n}\nonumber\\
&  \ \ \ \ \ \ \ \ \ \ \left(  \text{here, we renamed the summation index
}k\text{ as }q\right) \nonumber\\
&  =\left(  \sum_{q=1}^{n}\left(  -1\right)  ^{q+j}a_{i,q}\det\left(  A_{\sim
j,\sim q}\right)  \right)  _{1\leq i\leq n,\ 1\leq j\leq n}.
\label{pf.thm.adj.inverse.short.L}%
\end{align}

Now, we claim that%
\begin{equation}
\sum_{q=1}^{n}\left(  -1\right)  ^{q+j}a_{i,q}\det\left(  A_{\sim j,\sim
q}\right)  =\det A\cdot\delta_{i,j} \label{pf.thm.adj.inverse.short.twise}%
\end{equation}
for any $\left(  i,j\right)  \in\left\{  1,2,\ldots,n\right\}  ^{2}$.

[\textit{Proof of (\ref{pf.thm.adj.inverse.short.twise}):} Fix $\left(
i,j\right)  \in\left\{  1,2,\ldots,n\right\}  ^{2}$. We are in one of the
following two cases:

\textit{Case 1:} We have $i=j$.

\textit{Case 2:} We have $i\neq j$.

Let us consider Case 1 first. In this case, we have $i=j$. Hence,
$\delta_{i,j}=1$. Now, Theorem \ref{thm.laplace.gen} \textbf{(a)} (applied to
$p=i$) yields%
\[
\det A=\sum_{q=1}^{n}\underbrace{\left(  -1\right)  ^{i+q}}%
_{\substack{=\left(  -1\right)  ^{q+i}=\left(  -1\right)  ^{q+j}\\\text{(since
}i=j\text{)}}}a_{i,q}\det\left(  \underbrace{A_{\sim i,\sim q}}%
_{\substack{=A_{\sim j,\sim q}\\\text{(since }i=j\text{)}}}\right)
=\sum_{q=1}^{n}\left(  -1\right)  ^{q+j}a_{i,q}\det\left(  A_{\sim j,\sim
q}\right)  .
\]
In view of $\det A\cdot\underbrace{\delta_{i,j}}_{=1}=\det A$, this rewrites
as
\[
\det A\cdot\delta_{i,j}=\sum_{q=1}^{n}\left(  -1\right)  ^{q+j}a_{i,q}%
\det\left(  A_{\sim j,\sim q}\right)  .
\]
Thus, (\ref{pf.thm.adj.inverse.short.twise}) is proven in Case 1.

Let us next consider Case 2. In this case, we have $i\neq j$. Hence,
$\delta_{i,j}=0$ and $j\neq i$. Now, Proposition \ref{prop.laplace.0}
\textbf{(a)} (applied to $p=j$ and $r=i$) yields%
\[
0=\sum_{q=1}^{n}\underbrace{\left(  -1\right)  ^{j+q}}_{=\left(  -1\right)
^{q+j}}a_{i,q}\det\left(  A_{\sim j,\sim q}\right)  =\sum_{q=1}^{n}\left(
-1\right)  ^{q+j}a_{i,q}\det\left(  A_{\sim j,\sim q}\right)  .
\]
In view of $\det A\cdot\underbrace{\delta_{i,j}}_{=0}=0$, this rewrites as%
\[
\det A\cdot\delta_{i,j}=\sum_{q=1}^{n}\left(  -1\right)  ^{q+j}a_{i,q}%
\det\left(  A_{\sim j,\sim q}\right)  .
\]
Thus, (\ref{pf.thm.adj.inverse.short.twise}) is proven in Case 2.

We have now proven (\ref{pf.thm.adj.inverse.short.twise}) in each of the two
Cases 1 and 2. Thus, (\ref{pf.thm.adj.inverse.short.twise}) is proven.]

Now, (\ref{pf.thm.adj.inverse.short.L}) becomes%
\begin{align}
A\cdot\operatorname*{adj}A  &  =\left(  \underbrace{\sum_{q=1}^{n}\left(
-1\right)  ^{q+j}a_{i,q}\det\left(  A_{\sim j,\sim q}\right)  }%
_{\substack{=\det A\cdot\delta_{i,j}\\\text{(by
(\ref{pf.thm.adj.inverse.short.twise}))}}}\right)  _{1\leq i\leq n,\ 1\leq
j\leq n}\nonumber\\
&  =\left(  \det A\cdot\delta_{i,j}\right)  _{1\leq i\leq n,\ 1\leq j\leq
n}=\det A\cdot I_{n} \label{pf.thm.adj.inverse.short.part1}%
\end{align}
(by (\ref{pf.thm.adj.inverse.short.R})).

It now remains to prove that $\operatorname*{adj}A\cdot A=\det A\cdot I_{n}$.
One way to do this is by mimicking the above proof using Theorem
\ref{thm.laplace.gen} \textbf{(b)} and Proposition \ref{prop.laplace.0}
\textbf{(b)} instead of Theorem \ref{thm.laplace.gen} \textbf{(a)} and
Proposition \ref{prop.laplace.0} \textbf{(a)}. However, here is a slicker proof:

Let us forget that we fixed $A$. We thus have shown that
(\ref{pf.thm.adj.inverse.short.part1}) holds for every $n\times n$-matrix $A$.

Now, let $A$ be any $n\times n$-matrix. Then, we can apply
(\ref{pf.thm.adj.inverse.short.part1}) to $A^{T}$ instead of $A$. We thus
obtain%
\begin{equation}
A^{T}\cdot\operatorname*{adj}\left(  A^{T}\right)  =\underbrace{\det\left(
A^{T}\right)  }_{\substack{=\det A\\\text{(by Exercise \ref{exe.ps4.4})}%
}}\cdot I_{n}=\det A\cdot I_{n}. \label{pf.thm.adj.inverse.short.4}%
\end{equation}

Now, (\ref{pf.thm.adj.inverse.tranposes1}) (applied to $u=n$, $v=n$, $w=n$,
$P=\operatorname*{adj}A$ and $Q=A$) shows that%
\[
\left(  \operatorname*{adj}A\cdot A\right)  ^{T}=A^{T}\cdot\underbrace{\left(
\operatorname*{adj}A\right)  ^{T}}_{\substack{=\operatorname*{adj}\left(
A^{T}\right)  \\\text{(by Proposition \ref{prop.adj.transpose})}}}=A^{T}%
\cdot\operatorname*{adj}\left(  A^{T}\right)  =\det A\cdot I_{n}%
\ \ \ \ \ \ \ \ \ \ \left(  \text{by (\ref{pf.thm.adj.inverse.short.4}%
)}\right)  .
\]
Hence,%
\begin{align*}
\left(  \underbrace{\left(  \operatorname*{adj}A\cdot A\right)  ^{T}}_{=\det
A\cdot I_{n}}\right)  ^{T}  &  =\left(  \det A\cdot I_{n}\right)  ^{T}=\det
A\cdot\underbrace{\left(  I_{n}\right)  ^{T}}_{\substack{=I_{n}\\\text{(by
(\ref{pf.thm.adj.inverse.tranposes2}), applied to }u=n\text{)}}}\\
&  \ \ \ \ \ \ \ \ \ \ \left(  \text{by (\ref{pf.thm.adj.inverse.tranposes3}),
applied to }u=n\text{, }v=n\text{, }P=I_{n}\text{ and }\lambda=\det A\right)
\\
&  =\det A\cdot I_{n}.
\end{align*}
Compared with
\[
\left(  \left(  \operatorname*{adj}A\cdot A\right)  ^{T}\right)
^{T}=\operatorname*{adj}A\cdot A\ \ \ \ \ \ \ \ \ \ \left(  \text{by
(\ref{pf.thm.adj.inverse.tranposes4}), applied to }u=n\text{, }v=n\text{ and
}P=\operatorname*{adj}A\cdot A\right)  ,
\]
this yields $\operatorname*{adj}A\cdot A=\det A\cdot I_{n}$. Combined with
(\ref{pf.thm.adj.inverse.short.part1}), this yields%
\[
A\cdot\operatorname*{adj}A=\operatorname*{adj}A\cdot A=\det A\cdot I_{n}.
\]
This proves Theorem \ref{thm.adj.inverse}.
\end{proof}
\end{vershort}

\begin{verlong}
\begin{proof}
[Proof of Theorem \ref{thm.adj.inverse}.]For any two objects $i$ and $j$, we
define $\delta_{i,j}$ to be the element $%
\begin{cases}
1, & \text{if }i=j;\\
0, & \text{if }i\neq j
\end{cases}
$ of $\mathbb{K}$. Then, $I_{n}=\left(  \delta_{i,j}\right)  _{1\leq i\leq
n,\ 1\leq j\leq n}$ (by the definition of $I_{n}$), and thus%
\begin{equation}
\det A\cdot\underbrace{I_{n}}_{=\left(  \delta_{i,j}\right)  _{1\leq i\leq
n,\ 1\leq j\leq n}}=\det A\cdot\left(  \delta_{i,j}\right)  _{1\leq i\leq
n,\ 1\leq j\leq n}=\left(  \det A\cdot\delta_{i,j}\right)  _{1\leq i\leq
n,\ 1\leq j\leq n}. \label{pf.thm.adj.inverse.R}%
\end{equation}

On the other hand, let us write the matrix $A$ in the form $A=\left(
a_{i,j}\right)  _{1\leq i\leq n,\ 1\leq j\leq n}$. Then, the definition of the
product of two matrices shows that%
\begin{align}
A\cdot\operatorname*{adj}A  &  =\left(  \underbrace{\sum_{k=1}^{n}%
a_{i,k}\left(  -1\right)  ^{k+j}\det\left(  A_{\sim j,\sim k}\right)
}_{\substack{=\sum_{q=1}^{n}a_{i,q}\left(  -1\right)  ^{q+j}\det\left(
A_{\sim j,\sim q}\right)  \\\text{(here, we renamed the summation index
}k\text{ as }q\text{)}}}\right)  _{1\leq i\leq n,\ 1\leq j\leq n}\nonumber\\
&  \ \ \ \ \ \ \ \ \ \ \left(
\begin{array}
[c]{c}%
\text{since }A=\left(  a_{i,j}\right)  _{1\leq i\leq n,\ 1\leq j\leq n}\\
\text{and }\operatorname*{adj}A=\left(  \left(  -1\right)  ^{i+j}\det\left(
A_{\sim j,\sim i}\right)  \right)  _{1\leq i\leq n,\ 1\leq j\leq n}%
\end{array}
\right) \nonumber\\
&  =\left(  \sum_{q=1}^{n}\underbrace{a_{i,q}\left(  -1\right)  ^{q+j}%
}_{=\left(  -1\right)  ^{q+j}a_{i,q}}\det\left(  A_{\sim j,\sim q}\right)
\right)  _{1\leq i\leq n,\ 1\leq j\leq n}\nonumber\\
&  =\left(  \sum_{q=1}^{n}\left(  -1\right)  ^{q+j}a_{i,q}\det\left(  A_{\sim
j,\sim q}\right)  \right)  _{1\leq i\leq n,\ 1\leq j\leq n}.
\label{pf.thm.adj.inverse.L}%
\end{align}

Now, we claim that%
\begin{equation}
\sum_{q=1}^{n}\left(  -1\right)  ^{q+j}a_{i,q}\det\left(  A_{\sim j,\sim
q}\right)  =\det A\cdot\delta_{i,j} \label{pf.thm.adj.inverse.twise}%
\end{equation}
for any $\left(  i,j\right)  \in\left\{  1,2,\ldots,n\right\}  ^{2}$.

[\textit{Proof of (\ref{pf.thm.adj.inverse.twise}):} Fix $\left(  i,j\right)
\in\left\{  1,2,\ldots,n\right\}  ^{2}$. Thus, $i\in\left\{  1,2,\ldots
,n\right\}  $ and $j\in\left\{  1,2,\ldots,n\right\}  $. We are in one of the
following two cases:

\textit{Case 1:} We have $i=j$.

\textit{Case 2:} We have $i\neq j$.

Let us consider Case 1 first. In this case, we have $i=j$. Hence,
$\delta_{i,j}=1$ (by the definition of $\delta_{i,j}$). Now, Theorem
\ref{thm.laplace.gen} \textbf{(a)} (applied to $p=i$) yields%
\[
\det A=\sum_{q=1}^{n}\underbrace{\left(  -1\right)  ^{i+q}}%
_{\substack{=\left(  -1\right)  ^{q+i}=\left(  -1\right)  ^{q+j}\\\text{(since
}i=j\text{)}}}a_{i,q}\det\left(  \underbrace{A_{\sim i,\sim q}}%
_{\substack{=A_{\sim j,\sim q}\\\text{(since }i=j\text{)}}}\right)
=\sum_{q=1}^{n}\left(  -1\right)  ^{q+j}a_{i,q}\det\left(  A_{\sim j,\sim
q}\right)  .
\]
Hence,%
\[
\sum_{q=1}^{n}\left(  -1\right)  ^{q+j}a_{i,q}\det\left(  A_{\sim j,\sim
q}\right)  =\det A=\det A\cdot\delta_{i,j}%
\]
(since $\det A\cdot\underbrace{\delta_{i,j}}_{=1}=\det A$). Thus,
(\ref{pf.thm.adj.inverse.twise}) is proven in Case 1.

Let us next consider Case 2. In this case, we have $i\neq j$. Hence,
$\delta_{i,j}=0$ (by the definition of $\delta_{i,j}$) and $j\neq i$. Now,
Proposition \ref{prop.laplace.0} \textbf{(a)} (applied to $p=j$ and $r=i$)
yields%
\[
0=\sum_{q=1}^{n}\underbrace{\left(  -1\right)  ^{j+q}}_{=\left(  -1\right)
^{q+j}}a_{i,q}\det\left(  A_{\sim j,\sim q}\right)  =\sum_{q=1}^{n}\left(
-1\right)  ^{q+j}a_{i,q}\det\left(  A_{\sim j,\sim q}\right)  .
\]
Hence,%
\[
\sum_{q=1}^{n}\left(  -1\right)  ^{q+j}a_{i,q}\det\left(  A_{\sim j,\sim
q}\right)  =0=\det A\cdot\delta_{i,j}%
\]
(since $\det A\cdot\underbrace{\delta_{i,j}}_{=0}=0$). Thus,
(\ref{pf.thm.adj.inverse.twise}) is proven in Case 2.

We have now proven (\ref{pf.thm.adj.inverse.twise}) in each of the two Cases 1
and 2. Thus, (\ref{pf.thm.adj.inverse.twise}) always holds. This completes the
proof of (\ref{pf.thm.adj.inverse.twise}).]

Now, (\ref{pf.thm.adj.inverse.L}) becomes%
\begin{align}
A\cdot\operatorname*{adj}A  &  =\left(  \underbrace{\sum_{q=1}^{n}\left(
-1\right)  ^{q+j}a_{i,q}\det\left(  A_{\sim j,\sim q}\right)  }%
_{\substack{=\det A\cdot\delta_{i,j}\\\text{(by
(\ref{pf.thm.adj.inverse.twise}))}}}\right)  _{1\leq i\leq n,\ 1\leq j\leq
n}\nonumber\\
&  =\left(  \det A\cdot\delta_{i,j}\right)  _{1\leq i\leq n,\ 1\leq j\leq
n}=\det A\cdot I_{n} \label{pf.thm.adj.inverse.part1}%
\end{align}
(by (\ref{pf.thm.adj.inverse.R})).

It now remains to prove that $\operatorname*{adj}A\cdot A=\det A\cdot I_{n}$.
One way to do this is by mimicking the above proof using Theorem
\ref{thm.laplace.gen} \textbf{(b)} and Proposition \ref{prop.laplace.0}
\textbf{(b)} instead of Theorem \ref{thm.laplace.gen} \textbf{(a)} and
Proposition \ref{prop.laplace.0} \textbf{(a)}. However, here is a slicker proof:

Let us forget that we fixed $A$. We thus have shown that
(\ref{pf.thm.adj.inverse.part1}) holds for every $n\times n$-matrix $A$.

Now, let $A$ be any $n\times n$-matrix. Then, we can apply
(\ref{pf.thm.adj.inverse.part1}) to $A^{T}$ instead of $A$. We thus obtain%
\begin{equation}
A^{T}\cdot\operatorname*{adj}\left(  A^{T}\right)  =\underbrace{\det\left(
A^{T}\right)  }_{\substack{=\det A\\\text{(by Exercise \ref{exe.ps4.4})}%
}}\cdot I_{n}=\det A\cdot I_{n}. \label{pf.thm.adj.inverse.4}%
\end{equation}

Now, (\ref{pf.thm.adj.inverse.tranposes1}) (applied to $u=n$, $v=n$, $w=n$,
$P=\operatorname*{adj}A$ and $Q=A$) shows that%
\[
\left(  \operatorname*{adj}A\cdot A\right)  ^{T}=A^{T}\cdot\underbrace{\left(
\operatorname*{adj}A\right)  ^{T}}_{\substack{=\operatorname*{adj}\left(
A^{T}\right)  \\\text{(by Proposition \ref{prop.adj.transpose})}}}=A^{T}%
\cdot\operatorname*{adj}\left(  A^{T}\right)  =\det A\cdot I_{n}%
\ \ \ \ \ \ \ \ \ \ \left(  \text{by (\ref{pf.thm.adj.inverse.4})}\right)  .
\]
Hence,%
\begin{align*}
\left(  \underbrace{\left(  \operatorname*{adj}A\cdot A\right)  ^{T}}_{=\det
A\cdot I_{n}}\right)  ^{T}  &  =\left(  \det A\cdot I_{n}\right)  ^{T}=\det
A\cdot\underbrace{\left(  I_{n}\right)  ^{T}}_{\substack{=I_{n}\\\text{(by
(\ref{pf.thm.adj.inverse.tranposes2}), applied to }u=n\text{)}}}\\
&  \ \ \ \ \ \ \ \ \ \ \left(  \text{by (\ref{pf.thm.adj.inverse.tranposes3}),
applied to }u=n\text{, }v=n\text{, }P=I_{n}\text{ and }\lambda=\det A\right)
\\
&  =\det A\cdot I_{n}.
\end{align*}
Compared with
\[
\left(  \left(  \operatorname*{adj}A\cdot A\right)  ^{T}\right)
^{T}=\operatorname*{adj}A\cdot A\ \ \ \ \ \ \ \ \ \ \left(  \text{by
(\ref{pf.thm.adj.inverse.tranposes4}), applied to }u=n\text{, }v=n\text{ and
}P=\operatorname*{adj}A\cdot A\right)  ,
\]
this yields $\operatorname*{adj}A\cdot A=\det A\cdot I_{n}$. Combined with
(\ref{pf.thm.adj.inverse.part1}), this yields%
\[
A\cdot\operatorname*{adj}A=\operatorname*{adj}A\cdot A=\det A\cdot I_{n}.
\]
This proves Theorem \ref{thm.adj.inverse}.
\end{proof}
\end{verlong}

The following is a simple consequence of Theorem \ref{thm.adj.inverse}:

\begin{corollary}
\label{cor.adj.kernel}Let $n\in\mathbb{N}$. Let $A$ be an $n\times n$-matrix.
Let $v$ be a column vector with $n$ entries. If $Av=0_{n\times1}$, then $\det
A\cdot v=0_{n\times1}$.

(Recall that $0_{n\times1}$ denotes the $n\times1$ zero matrix, i.e., the
column vector with $n$ entries whose all entries are $0$.)
\end{corollary}

\begin{proof}
[Proof of Corollary \ref{cor.adj.kernel}.]Assume that $Av=0_{n\times1}$. It is
easy to see that every $m\in\mathbb{N}$ and every $n\times m$-matrix $B$
satisfy $I_{n}B=B$. Applying this to $m=1$ and $B=v$, we obtain $I_{n}v=v$.

It is also easy to see that every $m\in\mathbb{N}$ and every $m\times
n$-matrix $B$ satisfy $B\cdot0_{n\times1}=0_{m\times1}$. Applying this to
$m=n$ and $B=\operatorname*{adj}A$, we obtain $\operatorname*{adj}%
A\cdot0_{n\times1}=0_{n\times1}$.

Now, Theorem \ref{thm.adj.inverse} yields $\operatorname*{adj}A\cdot A=\det
A\cdot I_{n}$. Hence,%
\[
\underbrace{\left(  \operatorname*{adj}A\cdot A\right)  }_{=\det A\cdot I_{n}%
}v=\left(  \det A\cdot I_{n}\right)  v=\det A\cdot\underbrace{\left(
I_{n}v\right)  }_{=v}=\det A\cdot v.
\]
Compared to%
\begin{align*}
\left(  \operatorname*{adj}A\cdot A\right)  v  &  =\operatorname*{adj}%
A\cdot\underbrace{\left(  Av\right)  }_{=0_{n\times1}}%
\ \ \ \ \ \ \ \ \ \ \left(  \text{since matrix multiplication is
associative}\right) \\
&  =\operatorname*{adj}A\cdot0_{n\times1}=0_{n\times1},
\end{align*}
this yields $\det A\cdot v=0_{n\times1}$. This proves Corollary
\ref{cor.adj.kernel}.
\end{proof}

\begin{exercise}
\label{exe.adj(AB)}Let $n\in\mathbb{N}$. Let $A$ and $B$ be two $n\times
n$-matrices. Prove that%
\[
\operatorname*{adj}\left(  AB\right)  =\operatorname*{adj}B\cdot
\operatorname*{adj}A.
\]

\end{exercise}

Let me end this section with another application of Proposition
\ref{prop.laplace.0}:

\begin{exercise}
\label{exe.vander-hook}Let $n\in\mathbb{N}$. For every $n$ elements
$y_{1},y_{2},\ldots,y_{n}$ of $\mathbb{K}$, we define an element $V\left(
y_{1},y_{2},\ldots,y_{n}\right)  $ of $\mathbb{K}$ by
\[
V\left(  y_{1},y_{2},\ldots,y_{n}\right)  =\prod_{1\leq i<j\leq n}\left(
y_{i}-y_{j}\right)  .
\]

Let $x_{1},x_{2},\ldots,x_{n}$ be $n$ elements of $\mathbb{K}$. Let
$t\in\mathbb{K}$. Prove that%
\begin{align*}
&  \sum_{k=1}^{n}x_{k}V\left(  x_{1},x_{2},\ldots,x_{k-1},x_{k}+t,x_{k+1}%
,x_{k+2},\ldots,x_{n}\right) \\
&  =\left(  \dbinom{n}{2}t+\sum_{k=1}^{n}x_{k}\right)  V\left(  x_{1}%
,x_{2},\ldots,x_{n}\right)  .
\end{align*}

[\textbf{Hint:} Use Theorem \ref{thm.vander-det}, Laplace expansion,
Proposition \ref{prop.laplace.0} and the binomial formula.]
\end{exercise}

Exercise \ref{exe.vander-hook} is part of \cite[\S 4.3, Exercise
10]{Fulton-Young}.

\subsection{Inverting matrices}

We now will study inverses of matrices. We begin with a definition:

\begin{definition}
\label{def.matrices.inverses}Let $n\in\mathbb{N}$ and $m\in\mathbb{N}$. Let
$A$ be an $n\times m$-matrix.

\textbf{(a)} A \textit{left inverse} of $A$ means an $m\times n$-matrix $L$
such that $LA=I_{m}$. We say that the matrix $A$ is \textit{left-invertible}
if and only if a left inverse of $A$ exists.

\textbf{(b)} A \textit{right inverse} of $A$ means an $m\times n$-matrix $R$
such that $AR=I_{n}$. We say that the matrix $A$ is \textit{right-invertible}
if and only if a right inverse of $A$ exists.

\textbf{(c)} An \textit{inverse} of $A$ (or \textit{two-sided inverse} of $A$)
means an $m\times n$-matrix $B$ such that $BA=I_{m}$ and $AB=I_{n}$. We say
that the matrix $A$ is \textit{invertible} if and only if an inverse of $A$ exists.

The notions \textquotedblleft left-invertible\textquotedblright,
\textquotedblleft right-invertible\textquotedblright\ and \textquotedblleft
invertible\textquotedblright\ depend on the ring $\mathbb{K}$. We shall
therefore speak of \textquotedblleft left-invertible over $\mathbb{K}%
$\textquotedblright, \textquotedblleft right-invertible over $\mathbb{K}%
$\textquotedblright\ and \textquotedblleft invertible over $\mathbb{K}%
$\textquotedblright\ whenever the context does not unambiguously determine
$\mathbb{K}$.
\end{definition}

The notions of \textquotedblleft left inverse\textquotedblright,
\textquotedblleft right inverse\textquotedblright\ and \textquotedblleft
inverse\textquotedblright\ are not interchangeable (unlike for elements in a
commutative ring). We shall soon see in what cases they are identical; but
first, let us give a few examples.

\begin{example}
\label{exa.matrices.inverses}For this example, set $\mathbb{K}=\mathbb{Z}$.

Let $P$ be the $1\times2$-matrix $\left(
\begin{array}
[c]{cc}%
1 & 2
\end{array}
\right)  $. The matrix $P$ is right-invertible. For instance, $\left(
\begin{array}
[c]{c}%
-1\\
1
\end{array}
\right)  $ and $\left(
\begin{array}
[c]{c}%
3\\
-1
\end{array}
\right)  $ are two right inverses of $P$ (because $P\left(
\begin{array}
[c]{c}%
-1\\
1
\end{array}
\right)  =\left(
\begin{array}
[c]{c}%
1
\end{array}
\right)  =I_{1}$ and $P\left(
\begin{array}
[c]{c}%
3\\
-1
\end{array}
\right)  =\left(
\begin{array}
[c]{c}%
1
\end{array}
\right)  =I_{1}$). This example shows that the right inverse of a matrix is
not always unique.

The $2\times1$-matrix $P^{T}=\left(
\begin{array}
[c]{c}%
1\\
2
\end{array}
\right)  $ is left-invertible. The left inverses of $P^{T}$ are the transposes
of the right inverses of $P$.

The matrix $P$ is not left-invertible; the matrix $P^{T}$ is not right-invertible.

Let $Q$ be the $2\times2$-matrix $\left(
\begin{array}
[c]{cc}%
1 & -1\\
3 & -2
\end{array}
\right)  $. The matrix $Q$ is invertible. Its inverse is $\left(
\begin{array}
[c]{cc}%
-2 & 1\\
-3 & 1
\end{array}
\right)  $ (since $\left(
\begin{array}
[c]{cc}%
-2 & 1\\
-3 & 1
\end{array}
\right)  Q=I_{2}$ and $Q\left(
\begin{array}
[c]{cc}%
-2 & 1\\
-3 & 1
\end{array}
\right)  =I_{2}$). It is not hard to see that this is its only inverse.

Let $R$ be the $2\times2$-matrix $\left(
\begin{array}
[c]{cc}%
1 & 2\\
2 & -1
\end{array}
\right)  $. It can be seen that this matrix is not invertible \textbf{as a
matrix over }$\mathbb{Z}$. On the other hand, if we consider it as a matrix
over $\mathbb{K}=\mathbb{Q}$ instead, then it is invertible, with inverse
$\left(
\begin{array}
[c]{cc}%
1/5 & 2/5\\
2/5 & -1/5
\end{array}
\right)  $.
\end{example}

Of course, any inverse of a matrix $A$ is automatically both a left inverse of
$A$ and a right inverse of $A$. Thus, an invertible matrix $A$ is
automatically both left-invertible and right-invertible.

The following simple fact is an analogue of Proposition
\ref{prop.rings.inverse-uni}:

\begin{proposition}
\label{prop.matrices.inverse-uni}Let $n\in\mathbb{N}$ and $m\in\mathbb{N}$.
Let $A$ be an $n\times m$-matrix. Let $L$ be a left inverse of $A$. Let $R$ be
a right inverse of $A$.

\textbf{(a)} We have $L=R$.

\textbf{(b)} The matrix $A$ is invertible, and $L=R$ is an inverse of $A$.
\end{proposition}

\begin{proof}
[Proof of Proposition \ref{prop.matrices.inverse-uni}.]We know that $L$ is a
left inverse of $A$. In other words, $L$ is an $m\times n$-matrix such that
$LA=I_{m}$ (by the definition of a \textquotedblleft left
inverse\textquotedblright).

We know that $R$ is a right inverse of $A$. In other words, $R$ is an $m\times
n$-matrix such that $AR=I_{n}$ (by the definition of a \textquotedblleft right
inverse\textquotedblright).

Now, recall that $I_{m}G=G$ for every $k\in\mathbb{N}$ and every $m\times
k$-matrix $G$. Applying this to $k=n$ and $G=R$, we obtain $I_{m}R=R$.

Also, recall that $GI_{n}=G$ for every $k\in\mathbb{N}$ and every $k\times
n$-matrix $G$. Applying this to $k=m$ and $G=L$, we obtain $LI_{n}=L$. Thus,
$L=L\underbrace{I_{n}}_{=AR}=\underbrace{LA}_{=I_{m}}R=I_{m}R=R$. This proves
Proposition \ref{prop.matrices.inverse-uni} \textbf{(a)}.

\textbf{(b)} We have $LA=I_{m}$ and $A\underbrace{L}_{=R}=AR=I_{n}$. Thus, $L$
is an $m\times n$-matrix such that $LA=I_{m}$ and $AL=I_{n}$. In other words,
$L$ is an inverse of $A$ (by the definition of an \textquotedblleft
inverse\textquotedblright). Thus, $L=R$ is an inverse of $A$ (since $L=R$).
This proves Proposition \ref{prop.matrices.inverse-uni} \textbf{(b)}.
\end{proof}

\begin{corollary}
\label{cor.matrices.inverse-uni.cor}Let $n\in\mathbb{N}$ and $m\in\mathbb{N}$.
Let $A$ be an $n\times m$-matrix.

\textbf{(a)} If $A$ is left-invertible and right-invertible, then $A$ is invertible.

\textbf{(b)} If $A$ is invertible, then there exists exactly one inverse of
$A$.
\end{corollary}

\begin{proof}
[Proof of Corollary \ref{cor.matrices.inverse-uni.cor}.]\textbf{(a)} Assume
that $A$ is left-invertible and right-invertible. Thus, $A$ has a left inverse
$L$ (since $A$ is left-invertible). Consider this $L$. Also, $A$ has a right
inverse $R$ (since $A$ is right-invertible). Consider this $R$. Proposition
\ref{prop.matrices.inverse-uni} \textbf{(b)} yields that the matrix $A$ is
invertible, and $L=R$ is an inverse of $A$. Corollary
\ref{cor.matrices.inverse-uni.cor} \textbf{(a)} is proven.

\textbf{(b)} Assume that $A$ is invertible. Let $B$ and $B^{\prime}$ be any
two inverses of $A$. Since $B$ is an inverse of $A$, we know that $B$ is an
$m\times n$-matrix such that $BA=I_{m}$ and $AB=I_{n}$ (by the definition of
an \textquotedblleft inverse\textquotedblright). Thus, in particular, $B$ is
an $m\times n$-matrix such that $BA=I_{m}$. In other words, $B$ is a left
inverse of $A$. Since $B^{\prime}$ is an inverse of $A$, we know that
$B^{\prime}$ is an $m\times n$-matrix such that $B^{\prime}A=I_{m}$ and
$AB^{\prime}=I_{n}$ (by the definition of an \textquotedblleft
inverse\textquotedblright). Thus, in particular, $B^{\prime}$ is an $m\times
n$-matrix such that $AB^{\prime}=I_{n}$. In other words, $B^{\prime}$ is a
right inverse of $A$. Now, Proposition \ref{prop.matrices.inverse-uni}
\textbf{(a)} (applied to $L=B$ and $R=B^{\prime}$) shows that $B=B^{\prime}$.

Let us now forget that we fixed $B$ and $B^{\prime}$. We thus have shown that
if $B$ and $B^{\prime}$ are two inverses of $A$, then $B=B^{\prime}$. In other
words, any two inverses of $A$ are equal. In other words, there exists at most
one inverse of $A$. Since we also know that there exists at least one inverse
of $A$ (since $A$ is invertible), we thus conclude that there exists exactly
one inverse of $A$. This proves Corollary \ref{cor.matrices.inverse-uni.cor}
\textbf{(b)}.
\end{proof}

\begin{definition}
Let $n\in\mathbb{N}$ and $m\in\mathbb{N}$. Let $A$ be an invertible $n\times
m$-matrix. Corollary \ref{cor.matrices.inverse-uni.cor} \textbf{(b)} shows
that there exists exactly one inverse of $A$. Thus, we can speak of
\textquotedblleft\textit{the inverse of }$A$\textquotedblright. We denote this
inverse by $A^{-1}$.
\end{definition}

In contrast to Definition \ref{def.rings.invertible}, we do \textbf{not}
define the notation $B/A$ for two matrices $B$ and $A$ for which $A$ is
invertible. In fact, the trouble with such a notation would be its ambiguity:
should it mean $BA^{-1}$ or $A^{-1}B$ ? (In general, $BA^{-1}$ and $A^{-1}B$
are not the same.) Some authors do write $B/A$ for the matrices $BA^{-1}$ and
$A^{-1}B$ when these matrices are equal; but we shall not have a reason to do so.

\begin{remark}
\label{rmk.matrices.inverses.AA-1}Let $n\in\mathbb{N}$ and $m\in\mathbb{N}$.
Let $A$ be an invertible $n\times m$-matrix. Then, the inverse $A^{-1}$ of $A$
is an $m\times n$-matrix and satisfies $AA^{-1}=I_{n}$ and $A^{-1}A=I_{m}$.
This follows from the definition of the inverse of $A$; we are just stating it
once again, because it will later be used without mention.
\end{remark}

Example \ref{exa.matrices.inverses} (and your experiences with a linear
algebra class, if you have taken one) suggest the conjecture that only square
matrices can be invertible. Indeed, this is \textbf{almost} true. There is a
stupid counterexample: If $\mathbb{K}$ is a trivial ring, then every matrix
over $\mathbb{K}$ is invertible\footnote{For example, the $1\times2$-matrix
$\left(
\begin{array}
[c]{cc}%
0_{\mathbb{K}} & 0_{\mathbb{K}}%
\end{array}
\right)  $ over a trivial ring $\mathbb{K}$ is invertible, having inverse
$\left(
\begin{array}
[c]{c}%
0_{\mathbb{K}}\\
0_{\mathbb{K}}%
\end{array}
\right)  $. If you don't believe me, just check that%
\begin{align*}
\left(
\begin{array}
[c]{c}%
0_{\mathbb{K}}\\
0_{\mathbb{K}}%
\end{array}
\right)  \left(
\begin{array}
[c]{cc}%
0_{\mathbb{K}} & 0_{\mathbb{K}}%
\end{array}
\right)   &  =\left(
\begin{array}
[c]{cc}%
0_{\mathbb{K}} & 0_{\mathbb{K}}\\
0_{\mathbb{K}} & 0_{\mathbb{K}}%
\end{array}
\right)  =\left(
\begin{array}
[c]{cc}%
1_{\mathbb{K}} & 0_{\mathbb{K}}\\
0_{\mathbb{K}} & 1_{\mathbb{K}}%
\end{array}
\right)  \ \ \ \ \ \ \ \ \ \ \left(  \text{since }0_{\mathbb{K}}%
=1_{\mathbb{K}}\right) \\
&  =I_{2}%
\end{align*}
and $\left(
\begin{array}
[c]{cc}%
0_{\mathbb{K}} & 0_{\mathbb{K}}%
\end{array}
\right)  \left(
\begin{array}
[c]{c}%
0_{\mathbb{K}}\\
0_{\mathbb{K}}%
\end{array}
\right)  =\left(
\begin{array}
[c]{c}%
0_{\mathbb{K}}%
\end{array}
\right)  =\left(
\begin{array}
[c]{c}%
1_{\mathbb{K}}%
\end{array}
\right)  =I_{1}$.}. It turns out that this is the only case where nonsquare
matrices can be invertible. Indeed, we have the following:

\begin{theorem}
\label{thm.matrices.inverses.oblong}Let $n\in\mathbb{N}$ and $m\in\mathbb{N}$.
Let $A$ be an $n\times m$-matrix.

\textbf{(a)} If $A$ is left-invertible and if $n<m$, then $\mathbb{K}$ is a
trivial ring.

\textbf{(b)} If $A$ is right-invertible and if $n>m$, then $\mathbb{K}$ is a
trivial ring.

\textbf{(c)} If $A$ is invertible and if $n\neq m$, then $\mathbb{K}$ is a
trivial ring.
\end{theorem}

\begin{proof}
[Proof of Theorem \ref{thm.matrices.inverses.oblong}.]\textbf{(a)} Assume that
$A$ is left-invertible, and that $n<m$.

The matrix $A$ has a left inverse $L$ (since it is left-invertible). Consider
this $L$.

We know that $L$ is a left inverse of $A$. In other words, $L$ is an $m\times
n$-matrix such that $LA=I_{m}$ (by the definition of a \textquotedblleft left
inverse\textquotedblright). But (\ref{eq.exam.cauchy-binet.0}) (applied to
$m$, $n$, $L$ and $A$ instead of $n$, $m$, $A$ and $B$) yields $\det\left(
LA\right)  =0$ (since $n<m$). Thus, $0=\det\left(  \underbrace{LA}_{=I_{m}%
}\right)  =\det\left(  I_{m}\right)  =1$. Of course, the $0$ and the $1$ in
this equality mean the elements $0_{\mathbb{K}}$ and $1_{\mathbb{K}}$ of
$\mathbb{K}$ (rather than the integers $0$ and $1$); thus, it rewrites as
$0_{\mathbb{K}}=1_{\mathbb{K}}$. In other words, $\mathbb{K}$ is a trivial
ring. This proves Theorem \ref{thm.matrices.inverses.oblong} \textbf{(a)}.

\textbf{(b)} Assume that $A$ is right-invertible, and that $n>m$.

The matrix $A$ has a right inverse $R$ (since it is right-invertible).
Consider this $R$.

We know that $R$ is a right inverse of $A$. In other words, $R$ is an $m\times
n$-matrix such that $AR=I_{n}$ (by the definition of a \textquotedblleft right
inverse\textquotedblright). But (\ref{eq.exam.cauchy-binet.0}) (applied to
$B=R$) yields $\det\left(  AR\right)  =0$ (since $m<n$). Thus, $0=\det\left(
\underbrace{AR}_{=I_{n}}\right)  =\det\left(  I_{n}\right)  =1$. Of course,
the $0$ and the $1$ in this equality mean the elements $0_{\mathbb{K}}$ and
$1_{\mathbb{K}}$ of $\mathbb{K}$ (rather than the integers $0$ and $1$); thus,
it rewrites as $0_{\mathbb{K}}=1_{\mathbb{K}}$. In other words, $\mathbb{K}$
is a trivial ring. This proves Theorem \ref{thm.matrices.inverses.oblong}
\textbf{(b)}.

\textbf{(c)} Assume that $A$ is invertible, and that $n\neq m$. Since $n\neq
m$, we must be in one of the following two cases:

\textit{Case 1:} We have $n<m$.

\textit{Case 2:} We have $n>m$.

Let us first consider Case 1. In this case, we have $n<m$. Now, $A$ is
invertible, and thus left-invertible (since every invertible matrix is
left-invertible). Hence, $\mathbb{K}$ is a trivial ring (according to Theorem
\ref{thm.matrices.inverses.oblong} \textbf{(a)}). Thus, Theorem
\ref{thm.matrices.inverses.oblong} \textbf{(c)} is proven in Case 1.

Let us now consider Case 2. In this case, we have $n>m$. Now, $A$ is
invertible, and thus right-invertible (since every invertible matrix is
right-invertible). Hence, $\mathbb{K}$ is a trivial ring (according to Theorem
\ref{thm.matrices.inverses.oblong} \textbf{(b)}). Thus, Theorem
\ref{thm.matrices.inverses.oblong} \textbf{(c)} is proven in Case 2.

We have thus proven Theorem \ref{thm.matrices.inverses.oblong} \textbf{(c)} in
both Cases 1 and 2. Thus, Theorem \ref{thm.matrices.inverses.oblong}
\textbf{(c)} always holds.
\end{proof}

Theorem \ref{thm.matrices.inverses.oblong} \textbf{(c)} says that the question
whether a matrix is invertible is only interesting for square matrices, unless
the ring $\mathbb{K}$ is given so inexplicitly that we do not know whether it
is trivial or not\footnote{This actually happens rather often in algebra! For
example, rings are often defined by \textquotedblleft generators and
relations\textquotedblright\ (such as \textquotedblleft the ring with
commuting generators $a,b,c$ subject to the relations $a^{2}+b^{2}=c^{2}$ and
$ab=c$\textquotedblright). Sometimes the relations force the ring to become
trivial (for instance, the ring with generator $a$ and relations $a=1$ and
$a^{2}=2$ is clearly the trivial ring, because in this ring we have
$2=a^{2}=1^{2}=1$). Often this is not clear a-priori, and theorems such as
Theorem \ref{thm.matrices.inverses.oblong} can be used to show this. The
triviality of a ring can be a nontrivial statement! (Richman makes this point
in \cite{Richman}.)}. Let us now study the invertibility of a square matrix.
Here, the determinant turns out to be highly useful:

\begin{theorem}
\label{thm.matrices.inverses.square}Let $n\in\mathbb{N}$. Let $A$ be an
$n\times n$-matrix.

\textbf{(a)} The matrix $A$ is invertible if and only if the element $\det A$
of $\mathbb{K}$ is invertible (in $\mathbb{K}$).

\textbf{(b)} If $\det A$ is invertible, then the inverse of $A$ is
$A^{-1}=\dfrac{1}{\det A}\cdot\operatorname*{adj}A$.
\end{theorem}

When $\mathbb{K}$ is a field, the invertible elements of $\mathbb{K}$ are
precisely the nonzero elements of $\mathbb{K}$. Thus, when $\mathbb{K}$ is a
field, the statement of Theorem \ref{thm.matrices.inverses.square}
\textbf{(a)} can be rewritten as \textquotedblleft The matrix $A$ is
invertible if and only if $\det A\neq0$\textquotedblright; this is a
cornerstone of linear algebra. But our statement of Theorem
\ref{thm.matrices.inverses.square} \textbf{(a)} works for an arbitrary
commutative ring $\mathbb{K}$. In particular, it works for $\mathbb{K}%
=\mathbb{Z}$. Here is a consequence:

\begin{corollary}
\label{cor.matrices.inverses.square.ZZ}Let $n\in\mathbb{N}$. Let
$A\in\mathbb{Z}^{n\times n}$ be an $n\times n$-matrix over $\mathbb{Z}$. Then,
the matrix $A$ is invertible if and only if $\det A\in\left\{  1,-1\right\}  $.
\end{corollary}

\begin{proof}
[Proof of Corollary \ref{cor.matrices.inverses.square.ZZ}.]If $g$ is an
integer, then $g$ is invertible (in $\mathbb{Z}$) if and only if $g\in\left\{
1,-1\right\}  $. In other words, for every integer $g$, we have the following
equivalence:%
\begin{equation}
\left(  g\text{ is invertible (in }\mathbb{Z}\text{)}\right)
\Longleftrightarrow\left(  g\in\left\{  1,-1\right\}  \right)  .
\label{pf.cor.matrices.inverses.square.ZZ.1}%
\end{equation}

Now, Theorem \ref{thm.matrices.inverses.square} \textbf{(a)} (applied to
$\mathbb{K}=\mathbb{Z}$) yields that the matrix $A$ is invertible if and only
if the element $\det A$ of $\mathbb{Z}$ is invertible (in $\mathbb{Z}$). Thus,
we have the following chain of equivalences:%
\begin{align*}
&  \left(  \text{the matrix }A\text{ is invertible}\right) \\
&  \Longleftrightarrow\ \left(  \det A\text{ is invertible (in }%
\mathbb{Z}\text{)}\right)  \ \Longleftrightarrow\ \left(  \det A\in\left\{
1,-1\right\}  \right) \\
&  \ \ \ \ \ \ \ \ \ \ \left(  \text{by
(\ref{pf.cor.matrices.inverses.square.ZZ.1}), applied to }g=\det A\right)  .
\end{align*}
This proves Corollary \ref{cor.matrices.inverses.square.ZZ}.
\end{proof}

Notice that Theorem \ref{thm.matrices.inverses.square} \textbf{(b)} yields an
explicit way to compute the inverse of a square matrix $A$ (provided that we
can compute determinants and the inverse of $\det A$). This is not the fastest
way (at least not when $\mathbb{K}$ is a field), but it is useful for various
theoretical purposes.

\begin{proof}
[Proof of Theorem \ref{thm.matrices.inverses.square}.]\textbf{(a)}
$\Longrightarrow:$\ \ \ \ \footnote{In case you don't know what the notation
\textquotedblleft$\Longrightarrow:$\textquotedblright\ here means:
\par
Theorem \ref{thm.matrices.inverses.square} \textbf{(a)} is an
\textquotedblleft if and only if\textquotedblright\ assertion. In other words,
it asserts that $\mathcal{U}\Longleftrightarrow\mathcal{V}$ for two statements
$\mathcal{U}$ and $\mathcal{V}$. (In our case, $\mathcal{U}$ is the statement
\textquotedblleft the matrix $A$ is invertible\textquotedblright, and
$\mathcal{V}$ is the statement \textquotedblleft the element $\det A$ of
$\mathbb{K}$ is invertible (in $\mathbb{K}$)\textquotedblright.) In order to
prove a statement of the form $\mathcal{U}\Longleftrightarrow\mathcal{V}$, it
is sufficient to prove the implications $\mathcal{U}\Longrightarrow
\mathcal{V}$ and $\mathcal{U}\Longleftarrow\mathcal{V}$. Usually, these two
implications are proven separately (although not always; for instance, in the
proof of Corollary \ref{cor.matrices.inverses.square.ZZ}, we have used a chain
of equivalences to prove $\mathcal{U}\Longleftrightarrow\mathcal{V}$
directly). When writing such a proof, one often uses the abbreviations
\textquotedblleft$\Longrightarrow:$\textquotedblright\ and \textquotedblleft%
$\Longleftarrow:$\textquotedblright\ for \textquotedblleft Here comes the
proof of the implication $\mathcal{U}\Longrightarrow\mathcal{V}$%
:\textquotedblright\ and \textquotedblleft Here comes the proof of the
implication $\mathcal{U}\Longleftarrow\mathcal{V}$:\textquotedblright,
respectively.} Assume that the matrix $A$ is invertible. In other words, an
inverse $B$ of $A$ exists. Consider such a $B$.

The matrix $B$ is an inverse of $A$. In other words, $B$ is an $n\times
n$-matrix such that $BA=I_{n}$ and $AB=I_{n}$ (by the definition of an
\textquotedblleft inverse\textquotedblright). Theorem \ref{thm.det(AB)} yields
$\det\left(  AB\right)  =\det A\cdot\det B$, so that $\det A\cdot\det
B=\det\left(  \underbrace{AB}_{=I_{n}}\right)  =\det\left(  I_{n}\right)  =1$.
Of course, we also have $\det B\cdot\det A=\det A\cdot\det B=1$. Thus, $\det
B$ is an inverse of $\det A$ in $\mathbb{K}$. Therefore, the element $\det A$
is invertible (in $\mathbb{K}$). This proves the $\Longrightarrow$ direction
of Theorem \ref{thm.matrices.inverses.square} \textbf{(a)}.

$\Longleftarrow:$ Assume that the element $\det A$ is invertible (in
$\mathbb{K}$). Thus, its inverse $\dfrac{1}{\det A}$ exists. Theorem
\ref{thm.adj.inverse} yields%
\[
A\cdot\operatorname*{adj}A=\operatorname*{adj}A\cdot A=\det A\cdot I_{n}.
\]
Now, define an $n\times n$-matrix $B$ by $B=\dfrac{1}{\det A}\cdot
\operatorname*{adj}A$. Then,%
\[
A\underbrace{B}_{=\dfrac{1}{\det A}\cdot\operatorname*{adj}A}=A\cdot\left(
\dfrac{1}{\det A}\cdot\operatorname*{adj}A\right)  =\dfrac{1}{\det A}%
\cdot\underbrace{A\cdot\operatorname*{adj}A}_{=\det A\cdot I_{n}%
}=\underbrace{\dfrac{1}{\det A}\cdot\det A}_{=1}\cdot I_{n}=I_{n}%
\]
and%
\[
\underbrace{B}_{=\dfrac{1}{\det A}\cdot\operatorname*{adj}A}A=\dfrac{1}{\det
A}\cdot\underbrace{\operatorname*{adj}A\cdot A}_{=\det A\cdot I_{n}%
}=\underbrace{\dfrac{1}{\det A}\cdot\det A}_{=1}\cdot I_{n}=I_{n}.
\]

Thus, $B$ is an $n\times n$-matrix such that $BA=I_{n}$ and $AB=I_{n}$. In
other words, $B$ is an inverse of $A$ (by the definition of an
\textquotedblleft inverse\textquotedblright). Thus, an inverse of $A$ exists;
in other words, the matrix $A$ is invertible. This proves the $\Longleftarrow$
direction of Theorem \ref{thm.matrices.inverses.square} \textbf{(a)}.

We have now proven both directions of Theorem
\ref{thm.matrices.inverses.square} \textbf{(a)}. Theorem
\ref{thm.matrices.inverses.square} \textbf{(a)} is thus proven.

\textbf{(b)} Assume that $\det A$ is invertible. Thus, its inverse $\dfrac
{1}{\det A}$ exists. We define an $n\times n$-matrix $B$ by $B=\dfrac{1}{\det
A}\cdot\operatorname*{adj}A$. Then, $B$ is an inverse of $A$%
\ \ \ \ \footnote{We have shown this in our proof of the $\Longleftarrow$
direction of Theorem \ref{thm.matrices.inverses.square} \textbf{(a)}.}. In
other words, $B$ is \textbf{the} inverse of $A$. In other words, $B=A^{-1}$.
Hence, $A^{-1}=B=\dfrac{1}{\det A}\cdot\operatorname*{adj}A$. This proves
Theorem \ref{thm.matrices.inverses.square} \textbf{(b)}.
\end{proof}

\begin{corollary}
\label{cor.matrices.inverse.AB}Let $n\in\mathbb{N}$. Let $A$ and $B$ be two
$n\times n$-matrices such that $AB=I_{n}$.

\textbf{(a)} We have $BA=I_{n}$.

\textbf{(b)} The matrix $A$ is invertible, and the matrix $B$ is the inverse
of $A$.
\end{corollary}

\begin{proof}
[Proof of Corollary \ref{cor.matrices.inverse.AB}.]Theorem \ref{thm.det(AB)}
yields $\det\left(  AB\right)  =\det A\cdot\det B$, so that $\det A\cdot\det
B=\det\left(  \underbrace{AB}_{=I_{n}}\right)  =\det\left(  I_{n}\right)  =1$.
Of course, we also have $\det B\cdot\det A=\det A\cdot\det B=1$. Thus, $\det
B$ is an inverse of $\det A$ in $\mathbb{K}$. Therefore, the element $\det A$
is invertible (in $\mathbb{K}$). Therefore, the matrix $A$ is invertible
(according to the $\Longleftarrow$ direction of Theorem
\ref{thm.matrices.inverses.square} \textbf{(b)}). Thus, the inverse of $A$
exists. Let $C$ be this inverse. Thus, $C$ is a left inverse of $A$ (since
every inverse of $A$ is a left inverse of $A$).

The matrix $B$ is an $n\times n$-matrix satisfying $AB=I_{n}$. In other words,
$B$ is a right inverse of $A$. On the other hand, $C$ is a left inverse of
$A$. Hence, Proposition \ref{prop.matrices.inverse-uni} \textbf{(a)} (applied
to $L=C$ and $R=B$) yields $C=B$. Hence, the matrix $B$ is the inverse of $A$
(since the matrix $C$ is the inverse of $A$). Thus, Corollary
\ref{cor.matrices.inverse.AB} \textbf{(b)} is proven.

Since $B$ is the inverse of $A$, we have $BA=I_{n}$ and $AB=I_{n}$ (by the
definition of an \textquotedblleft inverse\textquotedblright). This proves
Corollary \ref{cor.matrices.inverse.AB} \textbf{(a)}.
\end{proof}

\begin{remark}
Corollary \ref{cor.matrices.inverse.AB} is \textbf{not} obvious! Matrix
multiplication, in general, is not commutative (we have $AB\neq BA$ more often
than not), and there is no reason to expect that $AB=I_{n}$ implies $BA=I_{n}%
$. The fact that this is nevertheless true for square matrices took us quite
some work to prove (we needed, among other things, the notion of an adjugate).
This fact would \textbf{not} hold for rectangular matrices. Nor does it hold
for \textquotedblleft infinite square matrices\textquotedblright: Without
wanting to go into the details of how products of infinite matrices are
defined, I invite you to check that the two infinite matrices $A=\left(
\begin{array}
[c]{cccc}%
0 & 1 & 0 & \cdots\\
0 & 0 & 1 & \cdots\\
0 & 0 & 0 & \cdots\\
\vdots & \vdots & \vdots & \ddots
\end{array}
\right)  $ and $B=A^{T}=\left(
\begin{array}
[c]{cccc}%
0 & 0 & 0 & \cdots\\
1 & 0 & 0 & \cdots\\
0 & 1 & 0 & \cdots\\
\vdots & \vdots & \vdots & \ddots
\end{array}
\right)  $ satisfy $AB=I_{\infty}$ but $BA\neq I_{\infty}$. This makes
Corollary \ref{cor.matrices.inverse.AB} \textbf{(a)} all the more interesting.
\end{remark}

Here are some more exercises involving matrices in \textquotedblleft
block-matrix form\textquotedblright:

\begin{exercise}
\label{exe.block2x2.VB+WD}Let $n\in\mathbb{N}$ and $m\in\mathbb{N}$. Let
$A\in\mathbb{K}^{n\times n}$, $B\in\mathbb{K}^{n\times m}$, $C\in
\mathbb{K}^{m\times n}$ and $D\in\mathbb{K}^{m\times m}$. Furthermore, let
$W\in\mathbb{K}^{m\times m}$ and $V\in\mathbb{K}^{m\times n}$ be such that
$VA=-WC$. Prove that%
\[
\det W\cdot\det\left(
\begin{array}
[c]{cc}%
A & B\\
C & D
\end{array}
\right)  =\det A\cdot\det\left(  VB+WD\right)  .
\]

[\textbf{Hint:} Use Exercise \ref{exe.block2x2.mult} to simplify the product
$\left(
\begin{array}
[c]{cc}%
I_{n} & 0_{n\times m}\\
V & W
\end{array}
\right)  \left(
\begin{array}
[c]{cc}%
A & B\\
C & D
\end{array}
\right)  $; then, take determinants.]
\end{exercise}

Exercise \ref{exe.block2x2.VB+WD} can often be used to compute the determinant
of a matrix given in block-matrix form (i.e., determinants of the form
$\det\left(
\begin{array}
[c]{cc}%
A & B\\
C & D
\end{array}
\right)  $) by only computing determinants of smaller matrices (such as $W$,
$A$ and $VB+WD$). It falls short of providing a general method for computing
such determinants\footnote{Indeed, it only gives a formula for $\det
W\cdot\det\left(
\begin{array}
[c]{cc}%
A & B\\
C & D
\end{array}
\right)  $, not for $\det\left(
\begin{array}
[c]{cc}%
A & B\\
C & D
\end{array}
\right)  $. If $\det W$ is invertible, then it allows for computing
$\det\left(
\begin{array}
[c]{cc}%
A & B\\
C & D
\end{array}
\right)  $; but Exercise \ref{exe.block2x2.VB+WD} gives no hint on how to find
matrices $W$ and $V$ such that $\det W$ is invertible and such that $VA=-WC$.
(Actually, such matrices do not always exist!)}, but it is one of the most
general facts about them. The next two exercises are two special cases of
Exercise \ref{exe.block2x2.VB+WD}:

\begin{exercise}
\label{exe.block2x2.schur}Let $n\in\mathbb{N}$ and $m\in\mathbb{N}$. Let
$A\in\mathbb{K}^{n\times n}$, $B\in\mathbb{K}^{n\times m}$, $C\in
\mathbb{K}^{m\times n}$ and $D\in\mathbb{K}^{m\times m}$ be such that the
matrix $A$ is invertible. Prove that%
\[
\det\left(
\begin{array}
[c]{cc}%
A & B\\
C & D
\end{array}
\right)  =\det A\cdot\det\left(  D-CA^{-1}B\right)  .
\]

\end{exercise}

Exercise \ref{exe.block2x2.schur} is known as the \textit{Schur complement
formula} (or, at least, it is one of several formulas sharing this name); and
the matrix $D-CA^{-1}B$ appearing on its right hand side is known as the
\textit{Schur complement} of the block $A$ in the matrix $\left(
\begin{array}
[c]{cc}%
A & B\\
C & D
\end{array}
\right)  $.

\begin{exercise}
\label{exe.block2x2.jacobi}Let $n\in\mathbb{N}$ and $m\in\mathbb{N}$. Let
$A\in\mathbb{K}^{n\times n}$, $B\in\mathbb{K}^{n\times m}$, $C\in
\mathbb{K}^{m\times n}$ and $D\in\mathbb{K}^{m\times m}$. Let $A^{\prime}%
\in\mathbb{K}^{n\times n}$, $B^{\prime}\in\mathbb{K}^{n\times m}$, $C^{\prime
}\in\mathbb{K}^{m\times n}$ and $D^{\prime}\in\mathbb{K}^{m\times m}$. Assume
that the matrix $\left(
\begin{array}
[c]{cc}%
A & B\\
C & D
\end{array}
\right)  $ is invertible, and that its inverse is the matrix $\left(
\begin{array}
[c]{cc}%
A^{\prime} & B^{\prime}\\
C^{\prime} & D^{\prime}%
\end{array}
\right)  $. Prove that%
\[
\det A=\det\left(
\begin{array}
[c]{cc}%
A & B\\
C & D
\end{array}
\right)  \cdot\det\left(  D^{\prime}\right)  .
\]

\end{exercise}

Exercise \ref{exe.block2x2.jacobi} can be rewritten in the following more
handy form:

\begin{exercise}
\label{exe.block2x2.jacobi.rewr}We shall use the notations introduced in
Definition \ref{def.submatrix}.

Let $n\in\mathbb{N}$. Let $A\in\mathbb{K}^{n\times n}$ be an invertible
matrix. Let $k\in\left\{  0,1,\ldots,n\right\}  $. Prove that%
\[
\det\left(  \operatorname*{sub}\nolimits_{1,2,\ldots,k}^{1,2,\ldots
,k}A\right)  =\det A\cdot\det\left(  \operatorname*{sub}%
\nolimits_{k+1,k+2,\ldots,n}^{k+1,k+2,\ldots,n}\left(  A^{-1}\right)  \right)
.
\]

\end{exercise}

Exercise \ref{exe.block2x2.jacobi.rewr} is a particular case of the so-called
\textit{Jacobi complementary minor theorem} (Exercise
\ref{addexe.jacobi-complement} further below).

\subsection{\label{sect.noncommring}Noncommutative rings}

I think that here is a good place to introduce two other basic notions from
algebra: that of a noncommutative ring, and that of a group.

\begin{definition}
\label{def.ring}The notion of a \textit{noncommutative ring} is defined in the
same way as we have defined a commutative ring (in Definition
\ref{def.commring}), except that we no longer require the \textquotedblleft
Commutativity of multiplication\textquotedblright\ axiom.
\end{definition}

As I have already said, the word \textquotedblleft
noncommutative\textquotedblright\ (in \textquotedblleft noncommutative
ring\textquotedblright) does not mean that commutativity of multiplication has
to be false in this ring; it only means that commutativity of multiplication
is not required. Thus, every commutative ring is a noncommutative ring.
Therefore, each of the examples of a commutative ring given in Section
\ref{sect.commring} is also an example of a noncommutative ring. Of course, it
is more interesting to see some examples of noncommutative rings which
actually fail to obey commutativity of multiplication. Here are some of these examples:

\begin{itemize}
\item If $n\in\mathbb{N}$ and if $\mathbb{K}$ is a commutative ring, then the
set $\mathbb{K}^{n\times n}$ of matrices becomes a noncommutative ring (when
endowed with the addition and multiplication of matrices, with the zero
$0_{n\times n}$ and with the unity $I_{n}$). This is actually a commutative
ring when $\mathbb{K}$ is trivial or when $n\leq1$, but in all
\textquotedblleft interesting\textquotedblright\ cases it is not commutative.

\item If you have heard of \href{https://en.wikipedia.org/wiki/Quaternion}{the
quaternions}, you should realize that they form a noncommutative ring.

\item Given a commutative ring $\mathbb{K}$ and $n$ distinct symbols
$X_{1},X_{2},\ldots,X_{n}$, we can define a \textit{ring of polynomials in the
\textbf{noncommutative} variables} $X_{1},X_{2},\ldots,X_{n}$ over
$\mathbb{K}$. We do not want to go into the details of its definition at this
point, but let us just mention some examples of its elements: For instance,
the ring of polynomials in the noncommutative variables $X$ and $Y$ over
$\mathbb{Q}$ contains elements such as $1+\dfrac{2}{3}X$, $X^{2}+\dfrac{3}%
{2}Y-7XY+YX$, $2XY$, $2YX$ and $5X^{2}Y-6XYX+7Y^{2}X$ (and of course, the
elements $XY$ and $YX$ are not equal).

\item If $n\in\mathbb{N}$ and if $\mathbb{K}$ is a commutative ring, then the
set of all lower-triangular $n\times n$-matrices over $\mathbb{K}$ becomes a
noncommutative ring (with addition, multiplication, zero and unity defined in
the same way as in $\mathbb{K}^{n\times n}$). This is because the sum and the
product of any two lower-triangular $n\times n$-matrices over $\mathbb{K}$ are
again lower-triangular\footnote{Check this! (For the sum, it is clear, but for
the product, it is an instructive exercise.)}, and because the matrices
$0_{n\times n}$ and $I_{n}$ are lower-triangular.

\item In contrast, the set of all invertible $2\times2$-matrices over
$\mathbb{K}$ is \textbf{not} a noncommutative ring (for example, because the
sum of the two invertible matrices $I_{2}$ and $-I_{2}$ is not
invertible\footnote{unless the ring $\mathbb{K}$ is trivial}).

\item If $\mathbb{K}$ is a commutative ring, then the set of all $3\times
3$-matrices (over $\mathbb{K}$) of the form $\left(
\begin{array}
[c]{ccc}%
a & b & c\\
0 & d & 0\\
0 & 0 & f
\end{array}
\right)  $ (with $a,b,c,d,f\in\mathbb{K}$) is a noncommutative ring (again,
with the same addition, multiplication, zero and unity as for $\mathbb{K}%
^{n\times n}$).\ \ \ \ \footnote{To check this, one needs to prove that the
matrices $0_{3\times3}$ and $I_{3}$ have this form, and that the sum and the
product of any two matrices of this form is again a matrix of this form. All
of this is clear, except for the claim about the product. The latter claim
follows from the computation%
\[
\left(
\begin{array}
[c]{ccc}%
a & b & c\\
0 & d & 0\\
0 & 0 & f
\end{array}
\right)  \left(
\begin{array}
[c]{ccc}%
a^{\prime} & b^{\prime} & c^{\prime}\\
0 & d^{\prime} & 0\\
0 & 0 & f^{\prime}%
\end{array}
\right)  =\left(
\begin{array}
[c]{ccc}%
aa^{\prime} & bd^{\prime}+ab^{\prime} & cf^{\prime}+ac^{\prime}\\
0 & dd^{\prime} & 0\\
0 & 0 & ff^{\prime}%
\end{array}
\right)  .
\]
}

\item On the other hand, if $\mathbb{K}$ is a commutative ring, then the set
of all $3\times3$-matrices (over $\mathbb{K}$) of the form $\left(
\begin{array}
[c]{ccc}%
a & b & 0\\
0 & c & d\\
0 & 0 & f
\end{array}
\right)  $ (with $a,b,c,d,f\in\mathbb{K}$) is \textbf{not} a noncommutative
ring (unless $\mathbb{K}$ is trivial), because products of matrices in this
set are not always in this set\footnote{Indeed, $\left(
\begin{array}
[c]{ccc}%
a & b & 0\\
0 & c & d\\
0 & 0 & f
\end{array}
\right)  \left(
\begin{array}
[c]{ccc}%
a^{\prime} & b^{\prime} & 0\\
0 & c^{\prime} & d^{\prime}\\
0 & 0 & f^{\prime}%
\end{array}
\right)  =\left(
\begin{array}
[c]{ccc}%
aa^{\prime} & ab^{\prime}+bc^{\prime} & bd^{\prime}\\
0 & cc^{\prime} & cd^{\prime}+df^{\prime}\\
0 & 0 & ff^{\prime}%
\end{array}
\right)  $ can have $bd^{\prime}\neq0$.}.
\end{itemize}

For the rest of this section, we let $\mathbb{L}$ be a \textbf{noncommutative}
ring. What can we do with elements of $\mathbb{L}$ ? We can do some of the
things that we can do with a commutative ring, but not all of them. For
example, we can still define the sum $a_{1}+a_{2}+\cdots+a_{n}$ and the
product $a_{1}a_{2}\cdots a_{n}$ of $n$ elements of a noncommutative ring. But
we cannot arbitrarily reorder the factors of a product and expect to always
get the same result! (With a sum, we can do this.) We can still define $na$
for any $n\in\mathbb{Z}$ and $a\in\mathbb{L}$ (in the same way as we defined
$na$ for $n\in\mathbb{Z}$ and $a\in\mathbb{K}$ when $\mathbb{K}$ was a
commutative ring). We can still define $a^{n}$ for any $n\in\mathbb{N}$ and
$a\in\mathbb{L}$ (again, in the same fashion as for commutative rings). The
identities (\ref{eq.rings.-(a+b)}), (\ref{eq.rings.-(-a)}),
(\ref{eq.rings.-(ab)}), (\ref{eq.rings.-(na)}), (\ref{eq.rings.nab}),
(\ref{eq.rings.nma}), (\ref{eq.rings.0**n}) and (\ref{eq.rings.a**(n+m)})
still hold when the commutative ring $\mathbb{K}$ is replaced by the
noncommutative ring $\mathbb{L}$; but the identities (\ref{eq.rings.-(ab)**n})
and (\ref{eq.rings.(a+b)**n}) may not (although they \textbf{do} hold if we
additionally assume that $ab=ba$). Finite sums such as $\sum_{s\in S}a_{s}$
(where $S$ is a finite set, and $a_{s}\in\mathbb{L}$ for every $s\in S$) are
well-defined, but finite products such as $\prod_{s\in S}a_{s}$ are not
(unless we specify the order in which their factors are to be multiplied).

We can define matrices over $\mathbb{L}$ in the same way as we have defined
matrices over $\mathbb{K}$. We can even define the determinant of a square
matrix over $\mathbb{L}$ using the formula (\ref{eq.det.eq.1}); however, this
determinant lacks many of the important properties that determinants over
$\mathbb{K}$ have (for instance, it satisfies neither Exercise \ref{exe.ps4.4}
nor Theorem \ref{thm.det(AB)}), and is therefore usually not
studied.\footnote{Some algebraists have come up with subtler notions of
determinants for matrices over noncommutative rings. But I don't want to go in
that direction here.}

We define the notion of an \textit{inverse} of an element $a\in\mathbb{L}$; in
order to do so, we simply replace $\mathbb{K}$ by $\mathbb{L}$ in Definition
\ref{def.rings.inverse}. (Now it suddenly matters that we required both $ab=1$
and $ba=1$ in Definition \ref{def.rings.inverse}.) Proposition
\ref{prop.rings.inverse-uni} still holds (and its proof still works) when
$\mathbb{K}$ is replaced by $\mathbb{L}$.

We define the notion of an \textit{invertible element} of $\mathbb{L}$; in
order to do so, we simply replace $\mathbb{K}$ by $\mathbb{L}$ in Definition
\ref{def.rings.invertible} \textbf{(a)}. We cannot directly replace
$\mathbb{K}$ by $\mathbb{L}$ in Definition \ref{def.rings.invertible}
\textbf{(b)}, because for two invertible elements $a$ and $b$ of $\mathbb{L}$
we do not necessarily have $\left(  ab\right)  ^{-1}=a^{-1}b^{-1}$; but
something very similar holds (namely, $\left(  ab\right)  ^{-1}=b^{-1}a^{-1}%
$). Trying to generalize Definition \ref{def.rings.invertible} \textbf{(c)} to
noncommutative rings is rather hopeless: In general, we cannot bring a
\textquotedblleft noncommutative fraction\textquotedblright\ of the form
$ba^{-1}+dc^{-1}$ to a \textquotedblleft common denominator\textquotedblright.

\begin{example}
\label{exa.rings.invertible-matrices}Let $\mathbb{K}$ be a commutative ring.
Let $n\in\mathbb{N}$. As we know, $\mathbb{K}^{n\times n}$ is a noncommutative
ring. The invertible elements of this ring are exactly the invertible $n\times
n$-matrices. (To see this, just compare the definition of an invertible
element of $\mathbb{K}^{n\times n}$ with the definition of an invertible
$n\times n$-matrix. These definitions are clearly equivalent.)
\end{example}

\subsection{Groups, and the group of units}

Let me finally define the notion of a \textit{group}.

\begin{definition}
\label{def.group}A \textit{group} means a set $G$ endowed with

\begin{itemize}
\item a binary operation called \textquotedblleft
multiplication\textquotedblright\ (or \textquotedblleft
composition\textquotedblright, or just \textquotedblleft binary
operation\textquotedblright), and denoted by $\cdot$, and written infix, and

\item an element called $1_{G}$ (or $e_{G}$)
\end{itemize}

\noindent such that the following axioms are satisfied:

\begin{itemize}
\item \textit{Associativity:} We have $a\left(  bc\right)  =\left(  ab\right)
c$ for all $a\in G$, $b\in G$ and $c\in G$. Here and in the following, the
expression \textquotedblleft$ab$\textquotedblright\ is shorthand for
\textquotedblleft$a\cdot b$\textquotedblright\ (as is usual for products of numbers).

\item \textit{Neutrality of }$1$\textit{:} We have $a1_{G}=1_{G}a=a$ for all
$a\in G$.

\item \textit{Existence of inverses:} For every $a\in G$, there exists an
element $a^{\prime}\in G$ such that $aa^{\prime}=a^{\prime}a=1_{G}$. This
$a^{\prime}$ is commonly denoted by $a^{-1}$ and called the \textit{inverse}
of $a$. (It is easy to check that it is unique.)
\end{itemize}
\end{definition}

\begin{definition}
The element $1_{G}$ of a group $G$ is denoted the \textit{neutral element} (or
the \textit{identity}) of $G$.

The binary operation $\cdot$ in Definition \ref{def.group} is usually not
identical with the binary operation $\cdot$ on the set of integers, and is
denoted by $\cdot_{G}$ when confusion can arise.
\end{definition}

The definition of a group has similarities with that of a noncommutative ring.
Viewed from a distance, it may look as if a noncommutative ring would
\textquotedblleft consist\textquotedblright\ of two groups with the same
underlying set. This is not quite correct, though, because the multiplication
in a nontrivial ring does not satisfy the \textquotedblleft existence of
inverses\textquotedblright\ axiom. But it is true that there are two groups in
every noncommutative ring:

\begin{proposition}
\label{prop.ring.groups}Let $\mathbb{L}$ be a noncommutative ring.

\textbf{(a)} The set $\mathbb{L}$, endowed with the \textbf{addition}
$+_{\mathbb{L}}$ (as multiplication) and the element $0_{\mathbb{L}}$ (as
neutral element), is a group. This group is called the \textit{additive group}
of $\mathbb{L}$, and denoted by $\mathbb{L}^{+}$.

\textbf{(b)} Let $\mathbb{L}^{\times}$ denote the set of all invertible
elements of $\mathbb{L}$. Then, the product of two elements of $\mathbb{L}%
^{\times}$ again belongs to $\mathbb{L}^{\times}$. Thus, we can define a
binary operation $\cdot_{\mathbb{L}^{\times}}$ on the set $\mathbb{L}^{\times
}$ (written infix) by
\[
a\cdot_{\mathbb{L}^{\times}}b=ab\ \ \ \ \ \ \ \ \ \ \text{for all }%
a\in\mathbb{L}^{\times}\text{ and }b\in\mathbb{L}^{\times}.
\]

The set $\mathbb{L}^{\times}$, endowed with the multiplication $\cdot
_{\mathbb{L}^{\times}}$ (as multiplication) and the element $1_{\mathbb{L}}$
(as neutral element), is a group. This group is called the \textit{group of
units} of $\mathbb{L}$.
\end{proposition}

\begin{proof}
[Proof of Proposition \ref{prop.ring.groups}.]\textbf{(a)} The addition
$+_{\mathbb{L}}$ is clearly a binary operation on $\mathbb{L}$, and the
element $0_{\mathbb{L}}$ is clearly an element of $\mathbb{L}$. The three
axioms in Definition \ref{def.group} are clearly satisfied for the binary
operation $+_{\mathbb{L}}$ and the element $0_{\mathbb{L}}$%
\ \ \ \ \footnote{In fact, they boil down to the \textquotedblleft
associativity of addition\textquotedblright, \textquotedblleft neutrality of
$0$\textquotedblright\ and \textquotedblleft existence of additive
inverses\textquotedblright\ axioms in the definition of a noncommutative
ring.}. Therefore, the set $\mathbb{L}$, endowed with the addition
$+_{\mathbb{L}}$ (as multiplication) and the element $0_{\mathbb{L}}$ (as
neutral element), is a group. This proves Proposition \ref{prop.ring.groups}
\textbf{(a)}.

\textbf{(b)} If $a\in\mathbb{L}^{\times}$ and $b\in\mathbb{L}^{\times}$, then
$ab\in\mathbb{L}^{\times}$\ \ \ \ \footnote{\textit{Proof.} Let $a\in
\mathbb{L}^{\times}$ and $b\in\mathbb{L}^{\times}$. We have $a\in
\mathbb{L}^{\times}$; in other words, $a$ is an invertible element of
$\mathbb{L}$ (because $\mathbb{L}^{\times}$ is the set of all invertible
elements of $\mathbb{L}$). Thus, the inverse $a^{-1}$ of $a$ is well-defined.
Similarly, the inverse $b^{-1}$ of $b$ is well-defined. Now, since we have%
\[
\left(  b^{-1}a^{-1}\right)  \left(  ab\right)  =b^{-1}\underbrace{a^{-1}%
a}_{=1_{\mathbb{L}}}b=b^{-1}b=1_{\mathbb{L}}%
\]
and%
\[
\left(  ab\right)  \left(  b^{-1}a^{-1}\right)  =a\underbrace{bb^{-1}%
}_{=1_{\mathbb{L}}}a^{-1}=aa^{-1}=1_{\mathbb{L}},
\]
we see that the element $b^{-1}a^{-1}$ of $\mathbb{L}$ is an inverse of $ab$.
Thus, the element $ab$ has an inverse. In other words, $ab$ is invertible. In
other words, $ab\in\mathbb{L}^{\times}$ (since $\mathbb{L}^{\times}$ is the
set of all invertible elements of $\mathbb{L}$), qed.}. In other words, the
product of two elements of $\mathbb{L}^{\times}$ again belongs to
$\mathbb{L}^{\times}$. Thus, we can define a binary operation $\cdot
_{\mathbb{L}^{\times}}$ on the set $\mathbb{L}^{\times}$ (written infix) by
\[
a\cdot_{\mathbb{L}^{\times}}b=ab\ \ \ \ \ \ \ \ \ \ \text{for all }%
a\in\mathbb{L}^{\times}\text{ and }b\in\mathbb{L}^{\times}.
\]
Also, $1_{\mathbb{L}}$ is an invertible element of $\mathbb{L}$ (indeed, its
inverse is $1_{\mathbb{L}}$), and thus an element of $\mathbb{L}^{\times}$.

Now, we need to prove that the set $\mathbb{L}^{\times}$, endowed with the
multiplication $\cdot_{\mathbb{L}^{\times}}$ (as multiplication) and the
element $1_{\mathbb{L}}$ (as neutral element), is a group. In order to do so,
we need to check that the \textquotedblleft associativity\textquotedblright,
\textquotedblleft neutrality of $1$\textquotedblright\ and \textquotedblleft
existence of inverses\textquotedblright\ axioms are satisfied.

The \textquotedblleft associativity\textquotedblright\ axiom follows from the
\textquotedblleft associativity of multiplication\textquotedblright\ axiom in
the definition of a noncommutative ring. The \textquotedblleft neutrality of
$1$\textquotedblright\ axiom follows from the \textquotedblleft
unitality\textquotedblright\ axiom in the definition of a noncommutative ring.
It thus remains to prove that the \textquotedblleft existence of
inverses\textquotedblright\ axiom holds.

Thus, let $a\in\mathbb{L}^{\times}$. We need to show that there exists an
$a^{\prime}\in\mathbb{L}^{\times}$ such that $a\cdot_{\mathbb{L}^{\times}%
}a^{\prime}=a^{\prime}\cdot_{\mathbb{L}^{\times}}a=1_{\mathbb{L}}$ (since
$1_{\mathbb{L}}$ is the neutral element of $\mathbb{L}^{\times}$).

We know that $a$ is an invertible element of $\mathbb{L}$ (sin${}$ce
$a\in\mathbb{L}^{\times}$); it thus has an inverse $a^{-1}$. Now, $a$ itself
is an inverse of $a^{-1}$ (since $aa^{-1}=1_{\mathbb{L}}$ and $a^{-1}%
a=1_{\mathbb{L}}$), and thus the element $a^{-1}$ of $\mathbb{L}$ has an
inverse. In other words, $a^{-1}$ is invertible, so that $a^{-1}\in
\mathbb{L}^{\times}$. The definition of the operation $\cdot_{\mathbb{L}%
^{\times}}$ shows that $a\cdot_{\mathbb{L}^{\times}}a^{-1}=aa^{-1}%
=1_{\mathbb{L}}$ and that $a^{-1}\cdot_{\mathbb{L}^{\times}}a=a^{-1}%
a=1_{\mathbb{L}}$. Hence, there exists an $a^{\prime}\in\mathbb{L}^{\times}$
such that $a\cdot_{\mathbb{L}^{\times}}a^{\prime}=a^{\prime}\cdot
_{\mathbb{L}^{\times}}a=1_{\mathbb{L}}$ (namely, $a^{\prime}=a^{-1}$). Thus we
have proven that the \textquotedblleft existence of inverses\textquotedblright%
\ axiom holds. The proof of Proposition \ref{prop.ring.groups} \textbf{(b)} is
thus complete.
\end{proof}

We now have a plentitude of examples of groups: For every noncommutative ring
$\mathbb{L}$, we have the two groups $\mathbb{L}^{+}$ and $\mathbb{L}^{\times
}$ defined in Proposition \ref{prop.ring.groups}. Another example, for every
set $X$, is the symmetric group of $X$ (endowed with the composition of
permutations as multiplication, and the identity permutation
$\operatorname*{id}:X\rightarrow X$ as the neutral element). (Many other
examples can be found in textbooks on algebra, such as \cite{Artin} or
\cite{Goodman}. Groups also naturally appear in the analysis of puzzles such
as Rubik's cube; this is explained in various sources such as \cite{Mulhol16},
\cite{Bump02} and \cite{Joyner08}, which can also be read as introductions to groups.)

\begin{remark}
Throwing all notational ballast aside, we can restate Proposition
\ref{prop.ring.groups} \textbf{(b)} as follows: The set of all invertible
elements of a noncommutative ring $\mathbb{L}$ is a group (where the binary
operation is multiplication). We can apply this to the case where
$\mathbb{L}=\mathbb{K}^{n\times n}$ for a commutative ring $\mathbb{K}$ and an
integer $n\in\mathbb{N}$. Thus, we obtain that the set of all invertible
elements of $\mathbb{K}^{n\times n}$ is a group. Since we know that the
invertible elements of $\mathbb{K}^{n\times n}$ are exactly the invertible
$n\times n$-matrices (by Example \ref{exa.rings.invertible-matrices}), we thus
have shown that the set of all invertible $n\times n$-matrices is a group.
This group is commonly denoted by $\operatorname*{GL}\nolimits_{n}\left(
\mathbb{K}\right)  $; it is called the \textit{general linear group of degree
}$n$ over $\mathbb{K}$.
\end{remark}

\subsection{Cramer's rule}

Let us return to the classical properties of determinants. We have already
proven many, but here is one more: It is an application of determinants to
solving systems of linear equations.

\begin{theorem}
\label{thm.cramer}Let $n\in\mathbb{N}$. Let $A$ be an $n\times n$-matrix. Let
$b=\left(  b_{1},b_{2},\ldots,b_{n}\right)  ^{T}$ be a column vector with $n$
entries (that is, an $n\times1$-matrix).\footnotemark

For every $j\in\left\{  1,2,\ldots,n\right\}  $, let $A_{j}^{\#}$ be the
$n\times n$-matrix obtained from $A$ by replacing the $j$-th column of $A$
with the vector $b$.

\textbf{(a)} We have $A\cdot\left(  \det\left(  A_{1}^{\#}\right)
,\det\left(  A_{2}^{\#}\right)  ,\ldots,\det\left(  A_{n}^{\#}\right)
\right)  ^{T}=\det A\cdot b$.

\textbf{(b)} Assume that the matrix $A$ is invertible. Then,%
\[
A^{-1}b=\left(  \dfrac{\det\left(  A_{1}^{\#}\right)  }{\det A},\dfrac
{\det\left(  A_{2}^{\#}\right)  }{\det A},\ldots,\dfrac{\det\left(  A_{n}%
^{\#}\right)  }{\det A}\right)  ^{T}.
\]

\end{theorem}

\footnotetext{The reader should keep in mind that $\left(  b_{1},b_{2}%
,\ldots,b_{n}\right)  ^{T}$ is just a space-saving way to write $\left(
\begin{array}
[c]{c}%
b_{1}\\
b_{2}\\
\vdots\\
b_{n}%
\end{array}
\right)  $.}Theorem \ref{thm.cramer} (or either part of it) is known as
\textit{Cramer's rule}.

\begin{remark}
A system of $n$ linear equations in $n$ variables $x_{1},x_{2},\ldots,x_{n}$
can be written in the form $Ax=b$, where $A$ is a fixed $n\times n$-matrix and
$b$ is a column vector with $n$ entries (and where $x$ is the column vector
$\left(  x_{1},x_{2},\ldots,x_{n}\right)  ^{T}$ containing all the variables).
When the matrix $A$ is invertible, it thus has a unique solution: namely,
$x=A^{-1}b$ (just multiply the equation $Ax=b$ from the left with $A^{-1}$ to
see this), and this solution can be computed using Theorem \ref{thm.cramer}.
This looks nice, but isn't actually all that useful for solving systems of
linear equations: For one thing, this does not immediately help us solve
systems of fewer or more than $n$ equations in $n$ variables; and even in the
case of exactly $n$ equations, the matrix $A$ coming from a system of linear
equations will not always be invertible (and in the more interesting cases, it
will not be). For another thing, at least when $\mathbb{K}$ is a field, there
are faster ways to solve a system of linear equations than anything that
involves computing $n+1$ determinants of $n\times n$-matrices. Theorem
\ref{thm.cramer} nevertheless turns out to be useful in proofs.
\end{remark}

\begin{vershort}
\begin{proof}
[Proof of Theorem \ref{thm.cramer}.]\textbf{(a)} Fix $j\in\left\{
1,2,\ldots,n\right\}  $. Let $C=A_{j}^{\#}$. Thus, $C=A_{j}^{\#}$ is the
$n\times n$-matrix obtained from $A$ by replacing the $j$-th column of $A$
with the vector $b$. In particular, the $j$-th column of $C$ is the vector
$b$. In other words, we have
\begin{equation}
c_{p,j}=b_{p}\ \ \ \ \ \ \ \ \ \ \text{for every }p\in\left\{  1,2,\ldots
,n\right\}  . \label{pf.thm.cramer.short.cpj}%
\end{equation}

Furthermore, the matrix $C$ is equal to the matrix $A$ in all columns but its
$j$-th column (because it is obtained from $A$ by replacing the $j$-th column
of $A$ with the vector $b$). Thus, if we cross out the $j$-th columns in both
matrices $C$ and $A$, then these two matrices become equal. Consequently,%
\begin{equation}
C_{\sim p,\sim j}=A_{\sim p,\sim j}\ \ \ \ \ \ \ \ \ \ \text{for every }%
p\in\left\{  1,2,\ldots,n\right\}  \label{pf.thm.cramer.short.CA}%
\end{equation}
(because the matrices $C_{\sim p,\sim j}$ and $A_{\sim p,\sim j}$ are obtained
by crossing out the $p$-th row and the $j$-th column in the matrices $C$ and
$A$, respectively). Now,%
\begin{align}
\det\left(  \underbrace{A_{j}^{\#}}_{=C}\right)   &  =\det C=\sum_{p=1}%
^{n}\left(  -1\right)  ^{p+j}\underbrace{c_{p,j}}_{\substack{=b_{p}\\\text{(by
(\ref{pf.thm.cramer.short.cpj}))}}}\det\left(  \underbrace{C_{\sim p,\sim j}%
}_{\substack{=A_{\sim p,\sim j}\\\text{(by (\ref{pf.thm.cramer.short.CA}))}%
}}\right) \nonumber\\
&  \ \ \ \ \ \ \ \ \ \ \left(
\begin{array}
[c]{c}%
\text{by Theorem \ref{thm.laplace.gen} \textbf{(b)}, applied}\\
\text{to }C\text{, }c_{i,j}\text{ and }j\text{ instead of }A\text{, }%
a_{i,j}\text{ and }q
\end{array}
\right) \nonumber\\
&  =\sum_{p=1}^{n}\left(  -1\right)  ^{p+j}b_{p}\det\left(  A_{\sim p,\sim
j}\right)  . \label{pf.thm.cramer.short.detAshj}%
\end{align}
Let us now forget that we fixed $j$. We thus have proven
(\ref{pf.thm.cramer.short.detAshj}) for every $j\in\left\{  1,2,\ldots
,n\right\}  $. Now, fix $i\in\left\{  1,2,\ldots,n\right\}  $. Then, for every
$p\in\left\{  1,2,\ldots,n\right\}  $ satisfying $p\neq i$, we have%
\begin{equation}
\sum_{q=1}^{n}b_{p}\left(  -1\right)  ^{p+q}a_{i,q}\det\left(  A_{\sim p,\sim
q}\right)  =0 \label{pf.thm.cramer.short.termkiller1}%
\end{equation}
\footnote{\textit{Proof of (\ref{pf.thm.cramer.short.termkiller1}):} Let
$p\in\left\{  1,2,\ldots,n\right\}  $ be such that $p\neq i$. Hence,
Proposition \ref{prop.laplace.0} \textbf{(a)} (applied to $r=i$) shows that%
\begin{equation}
0=\sum_{q=1}^{n}\left(  -1\right)  ^{p+q}a_{i,q}\det\left(  A_{\sim p,\sim
q}\right)  . \label{pf.thm.cramer.short.termkiller1.pf.1}%
\end{equation}
Now,%
\[
\sum_{q=1}^{n}b_{p}\left(  -1\right)  ^{p+q}a_{i,q}\det\left(  A_{\sim p,\sim
q}\right)  =b_{p}\underbrace{\sum_{q=1}^{n}\left(  -1\right)  ^{p+q}%
a_{i,q}\det\left(  A_{\sim p,\sim q}\right)  }_{\substack{=0\\\text{(by
(\ref{pf.thm.cramer.short.termkiller1.pf.1}))}}}=0.
\]
Thus, (\ref{pf.thm.cramer.short.termkiller1}) is proven.}. Also, we have%
\begin{equation}
\sum_{q=1}^{n}b_{i}\left(  -1\right)  ^{i+q}a_{i,q}\det\left(  A_{\sim i,\sim
q}\right)  =\det A\cdot b_{i} \label{pf.thm.cramer.short.termkiller2}%
\end{equation}
\footnote{\textit{Proof of (\ref{pf.thm.cramer.short.termkiller2}):} Applying
Theorem \ref{thm.laplace.gen} \textbf{(a)} to $p=i$, we obtain%
\begin{equation}
\det A=\sum_{q=1}^{n}\left(  -1\right)  ^{i+q}a_{i,q}\det\left(  A_{\sim
i,\sim q}\right)  . \label{pf.thm.cramer.short.termkiller2.pf.1}%
\end{equation}
Now,%
\[
\sum_{q=1}^{n}b_{i}\left(  -1\right)  ^{i+q}a_{i,q}\det\left(  A_{\sim i,\sim
q}\right)  =b_{i}\underbrace{\sum_{q=1}^{n}\left(  -1\right)  ^{i+q}%
a_{i,q}\det\left(  A_{\sim i,\sim q}\right)  }_{\substack{=\det A\\\text{(by
(\ref{pf.thm.cramer.short.termkiller2.pf.1}))}}}=b_{i}\det A=\det A\cdot
b_{i}.
\]
This proves (\ref{pf.thm.cramer.short.termkiller2}).}. Now,%
\begin{align}
&  \sum_{k=1}^{n}a_{i,k}\underbrace{\det\left(  A_{k}^{\#}\right)
}_{\substack{=\sum_{p=1}^{n}\left(  -1\right)  ^{p+k}b_{p}\det\left(  A_{\sim
p,\sim k}\right)  \\\text{(by (\ref{pf.thm.cramer.short.detAshj}), applied to
}j=k\text{)}}}\nonumber\\
&  =\sum_{k=1}^{n}a_{i,k}\sum_{p=1}^{n}\left(  -1\right)  ^{p+k}b_{p}%
\det\left(  A_{\sim p,\sim k}\right)  =\sum_{q=1}^{n}a_{i,q}\sum_{p=1}%
^{n}\left(  -1\right)  ^{p+q}b_{p}\det\left(  A_{\sim p,\sim q}\right)
\nonumber\\
&  \ \ \ \ \ \ \ \ \ \ \left(  \text{here, we renamed the summation index
}k\text{ as }q\text{ in the first sum}\right) \nonumber\\
&  =\underbrace{\sum_{q=1}^{n}\sum_{p=1}^{n}}_{\substack{=\sum_{p=1}^{n}%
\sum_{q=1}^{n}\\=\sum_{p\in\left\{  1,2,\ldots,n\right\}  }\sum_{q=1}^{n}%
}}\underbrace{a_{i,q}\left(  -1\right)  ^{p+q}b_{p}}_{=b_{p}\left(  -1\right)
^{p+q}a_{i,q}}\det\left(  A_{\sim p,\sim q}\right) \nonumber\\
&  =\sum_{p\in\left\{  1,2,\ldots,n\right\}  }\sum_{q=1}^{n}b_{p}\left(
-1\right)  ^{p+q}a_{i,q}\det\left(  A_{\sim p,\sim q}\right) \nonumber\\
&  =\sum_{\substack{p\in\left\{  1,2,\ldots,n\right\}  ;\\p\neq i}%
}\underbrace{\sum_{q=1}^{n}b_{p}\left(  -1\right)  ^{p+q}a_{i,q}\det\left(
A_{\sim p,\sim q}\right)  }_{\substack{=0\\\text{(by
(\ref{pf.thm.cramer.short.termkiller1}))}}}+\underbrace{\sum_{q=1}^{n}%
b_{i}\left(  -1\right)  ^{i+q}a_{i,q}\det\left(  A_{\sim i,\sim q}\right)
}_{\substack{=\det A\cdot b_{i}\\\text{(by
(\ref{pf.thm.cramer.short.termkiller2}))}}}\nonumber\\
&  \ \ \ \ \ \ \ \ \ \ \left(  \text{here, we have split off the addend for
}p=i\text{ from the sum}\right) \nonumber\\
&  =\underbrace{\sum_{\substack{p\in\left\{  1,2,\ldots,n\right\}  ;\\p\neq
i}}0}_{=0}+\det A\cdot b_{i}=\det A\cdot b_{i}.
\label{pf.thm.cramer.short.almostthere}%
\end{align}
Now, let us forget that we fixed $i$. We thus have proven
(\ref{pf.thm.cramer.short.almostthere}) for every $i\in\left\{  1,2,\ldots
,n\right\}  $. Now, let $d$ be the vector $\left(  \det\left(  A_{1}%
^{\#}\right)  ,\det\left(  A_{2}^{\#}\right)  ,\ldots,\det\left(  A_{n}%
^{\#}\right)  \right)  ^{T}$. Thus,%
\begin{align*}
d  &  =\left(  \det\left(  A_{1}^{\#}\right)  ,\det\left(  A_{2}^{\#}\right)
,\ldots,\det\left(  A_{n}^{\#}\right)  \right)  ^{T}=\left(
\begin{array}
[c]{c}%
\det\left(  A_{1}^{\#}\right) \\
\det\left(  A_{2}^{\#}\right) \\
\vdots\\
\det\left(  A_{n}^{\#}\right)
\end{array}
\right) \\
&  =\left(  \det\left(  A_{i}^{\#}\right)  \right)  _{1\leq i\leq n,\ 1\leq
j\leq1}.
\end{align*}
The definition of the product of two matrices shows that%
\begin{align*}
A\cdot d  &  =\left(  \underbrace{\sum_{k=1}^{n}a_{i,k}\det\left(  A_{k}%
^{\#}\right)  }_{\substack{=\det A\cdot b_{i}\\\text{(by
(\ref{pf.thm.cramer.short.almostthere}))}}}\right)  _{1\leq i\leq n,\ 1\leq
j\leq1}\\
&  \ \ \ \ \ \ \ \ \ \ \left(  \text{since }A=\left(  a_{i,j}\right)  _{1\leq
i\leq n,\ 1\leq j\leq n}\text{ and }d=\left(  \det\left(  A_{i}^{\#}\right)
\right)  _{1\leq i\leq n,\ 1\leq j\leq1}\right) \\
&  =\left(  \det A\cdot b_{i}\right)  _{1\leq i\leq n,\ 1\leq j\leq1}=\left(
\det A\cdot b_{1},\det A\cdot b_{2},\ldots,\det A\cdot b_{n}\right)  ^{T}.
\end{align*}
Comparing this with%
\[
\det A\cdot\underbrace{b}_{=\left(  b_{1},b_{2},\ldots,b_{n}\right)  ^{T}%
}=\det A\cdot\left(  b_{1},b_{2},\ldots,b_{n}\right)  ^{T}=\left(  \det A\cdot
b_{1},\det A\cdot b_{2},\ldots,\det A\cdot b_{n}\right)  ^{T},
\]
we obtain $A\cdot d=\det A\cdot b$. Since $d=\left(  \det\left(  A_{1}%
^{\#}\right)  ,\det\left(  A_{2}^{\#}\right)  ,\ldots,\det\left(  A_{n}%
^{\#}\right)  \right)  ^{T}$, we can rewrite this as $A\cdot\left(
\det\left(  A_{1}^{\#}\right)  ,\det\left(  A_{2}^{\#}\right)  ,\ldots
,\det\left(  A_{n}^{\#}\right)  \right)  ^{T}=\det A\cdot b$. This proves
Theorem \ref{thm.cramer} \textbf{(a)}.

\textbf{(b)} Theorem \ref{thm.matrices.inverses.square} \textbf{(a)} shows
that the matrix $A$ is invertible if and only if the element $\det A$ of
$\mathbb{K}$ is invertible (in $\mathbb{K}$). Hence, the element $\det A$ of
$\mathbb{K}$ is invertible (since the matrix $A$ is invertible). Thus,
$\dfrac{1}{\det A}$ is well-defined. Clearly, $\underbrace{\dfrac{1}{\det
A}\cdot\det A}_{=1}\cdot b=b$, so that%
\begin{align*}
b  &  =\dfrac{1}{\det A}\cdot\underbrace{\det A\cdot b}_{\substack{=A\cdot
\left(  \det\left(  A_{1}^{\#}\right)  ,\det\left(  A_{2}^{\#}\right)
,\ldots,\det\left(  A_{n}^{\#}\right)  \right)  ^{T}\\\text{(by Theorem
\ref{thm.cramer} \textbf{(a)})}}}\\
&  =\dfrac{1}{\det A}\cdot A\cdot\left(  \det\left(  A_{1}^{\#}\right)
,\det\left(  A_{2}^{\#}\right)  ,\ldots,\det\left(  A_{n}^{\#}\right)
\right)  ^{T}\\
&  =A\cdot\underbrace{\left(  \dfrac{1}{\det A}\cdot\left(  \det\left(
A_{1}^{\#}\right)  ,\det\left(  A_{2}^{\#}\right)  ,\ldots,\det\left(
A_{n}^{\#}\right)  \right)  ^{T}\right)  }_{\substack{=\left(  \dfrac{1}{\det
A}\det\left(  A_{1}^{\#}\right)  ,\dfrac{1}{\det A}\det\left(  A_{2}%
^{\#}\right)  ,\ldots,\dfrac{1}{\det A}\det\left(  A_{n}^{\#}\right)  \right)
^{T}\\=\left(  \dfrac{\det\left(  A_{1}^{\#}\right)  }{\det A},\dfrac
{\det\left(  A_{2}^{\#}\right)  }{\det A},\ldots,\dfrac{\det\left(  A_{n}%
^{\#}\right)  }{\det A}\right)  ^{T}}}\\
&  =A\cdot\left(  \dfrac{\det\left(  A_{1}^{\#}\right)  }{\det A},\dfrac
{\det\left(  A_{2}^{\#}\right)  }{\det A},\ldots,\dfrac{\det\left(  A_{n}%
^{\#}\right)  }{\det A}\right)  ^{T}.
\end{align*}
Therefore,
\begin{align*}
&  A^{-1}\underbrace{b}_{=A\cdot\left(  \dfrac{\det\left(  A_{1}^{\#}\right)
}{\det A},\dfrac{\det\left(  A_{2}^{\#}\right)  }{\det A},\ldots,\dfrac
{\det\left(  A_{n}^{\#}\right)  }{\det A}\right)  ^{T}}\\
&  =\underbrace{A^{-1}A}_{=I_{n}}\cdot\left(  \dfrac{\det\left(  A_{1}%
^{\#}\right)  }{\det A},\dfrac{\det\left(  A_{2}^{\#}\right)  }{\det A}%
,\ldots,\dfrac{\det\left(  A_{n}^{\#}\right)  }{\det A}\right)  ^{T}\\
&  =I_{n}\cdot\left(  \dfrac{\det\left(  A_{1}^{\#}\right)  }{\det A}%
,\dfrac{\det\left(  A_{2}^{\#}\right)  }{\det A},\ldots,\dfrac{\det\left(
A_{n}^{\#}\right)  }{\det A}\right)  ^{T}=\left(  \dfrac{\det\left(
A_{1}^{\#}\right)  }{\det A},\dfrac{\det\left(  A_{2}^{\#}\right)  }{\det
A},\ldots,\dfrac{\det\left(  A_{n}^{\#}\right)  }{\det A}\right)  ^{T}.
\end{align*}
This proves Theorem \ref{thm.cramer} \textbf{(b)}.
\end{proof}
\end{vershort}

\begin{verlong}
\begin{proof}
[Proof of Theorem \ref{thm.cramer}.]\textbf{(a)} Fix $j\in\left\{
1,2,\ldots,n\right\}  $. Let $C=A_{j}^{\#}$.

We know that $A_{j}^{\#}$ is the $n\times n$-matrix obtained from $A$ by
replacing the $j$-th column of $A$ with the vector $b$. In other words, $C$ is
the $n\times n$-matrix obtained from $A$ by replacing the $j$-th column of $A$
with the vector $b$ (because $C=A_{j}^{\#}$). In other words, we have%
\begin{align}
&  \left(  \left(  \text{the }u\text{-th column of }C\right)  =\left(
\text{the }u\text{-th column of }A\right)  \right.  \label{pf.thm.cramer.C.1}%
\\
&  \ \ \ \ \ \ \ \ \ \ \left.  \text{for all }u\in\left\{  1,2,\ldots
,n\right\}  \text{ satisfying }u\neq j\right)  ,\nonumber
\end{align}
whereas%
\begin{equation}
\left(  \text{the }j\text{-th column of }C\right)  =b.
\label{pf.thm.cramer.C.2}%
\end{equation}

Let us write the $n\times n$-matrix $A$ in the form $A=\left(  a_{i,j}\right)
_{1\leq i\leq n,\ 1\leq j\leq n}$. Thus, for every $u\in\left\{
1,2,\ldots,n\right\}  $, we have%
\begin{equation}
\left(  \text{the }u\text{-th column of }A\right)  =\left(  a_{1,u}%
,a_{2,u},\ldots,a_{n,u}\right)  ^{T}. \label{pf.thm.cramer.u-th-col-A}%
\end{equation}

On the other hand, let us write the $n\times n$-matrix $C$ in the form
$C=\left(  c_{i,j}\right)  _{1\leq i\leq n,\ 1\leq j\leq n}$. Thus, for every
$u\in\left\{  1,2,\ldots,n\right\}  $, we have%
\begin{equation}
\left(  \text{the }u\text{-th column of }C\right)  =\left(  c_{1,u}%
,c_{2,u},\ldots,c_{n,u}\right)  ^{T}. \label{pf.thm.cramer.u-th-col-C}%
\end{equation}

Now, it is easy to see that
\begin{equation}
c_{p,j}=b_{p}\ \ \ \ \ \ \ \ \ \ \text{for every }p\in\left\{  1,2,\ldots
,n\right\}  \label{pf.thm.cramer.cpj}%
\end{equation}
\footnote{\textit{Proof of (\ref{pf.thm.cramer.cpj}):} Applying
(\ref{pf.thm.cramer.u-th-col-C}) to $u=j$, we obtain
\[
\left(  \text{the }j\text{-th column of }C\right)  =\left(  c_{1,j}%
,c_{2,j},\ldots,c_{n,j}\right)  ^{T}.
\]
Thus,%
\begin{align*}
\left(  c_{1,j},c_{2,j},\ldots,c_{n,j}\right)  ^{T}  &  =\left(  \text{the
}j\text{-th column of }C\right)  =b\ \ \ \ \ \ \ \ \ \ \left(  \text{by
(\ref{pf.thm.cramer.C.2})}\right) \\
&  =\left(  b_{1},b_{2},\ldots,b_{n}\right)  ^{T}.
\end{align*}
Thus,%
\[
\left(  c_{1,j},c_{2,j},\ldots,c_{n,j}\right)  =\left(  \underbrace{\left(
c_{1,j},c_{2,j},\ldots,c_{n,j}\right)  ^{T}}_{=\left(  b_{1},b_{2}%
,\ldots,b_{n}\right)  ^{T}}\right)  ^{T}=\left(  \left(  b_{1},b_{2}%
,\ldots,b_{n}\right)  ^{T}\right)  ^{T}=\left(  b_{1},b_{2},\ldots
,b_{n}\right)  .
\]
In other words, $c_{p,j}=b_{p}$ for every $p\in\left\{  1,2,\ldots,n\right\}
$. This proves (\ref{pf.thm.cramer.cpj}).}. Furthermore,%
\begin{equation}
c_{p,q}=a_{p,q}\ \ \ \ \ \ \ \ \ \ \text{for every }p\in\left\{
1,2,\ldots,n\right\}  \text{ and }q\in\left\{  1,2,\ldots,n\right\}  \text{
satisfying }q\neq j \label{pf.thm.cramer.cpq}%
\end{equation}
\footnote{\textit{Proof of (\ref{pf.thm.cramer.cpq}):} Let $q\in\left\{
1,2,\ldots,n\right\}  $ be such that $q\neq j$. We need to prove that
$c_{p,q}=a_{p,q}$ for every $p\in\left\{  1,2,\ldots,n\right\}  $.
\par
Applying (\ref{pf.thm.cramer.u-th-col-C}) to $u=q$, we obtain
\[
\left(  \text{the }q\text{-th column of }C\right)  =\left(  c_{1,q}%
,c_{2,q},\ldots,c_{n,q}\right)  ^{T}.
\]
Thus,%
\begin{align*}
\left(  c_{1,q},c_{2,q},\ldots,c_{n,q}\right)  ^{T}  &  =\left(  \text{the
}q\text{-th column of }C\right) \\
&  =\left(  \text{the }q\text{-th column of }A\right)
\ \ \ \ \ \ \ \ \ \ \left(  \text{by (\ref{pf.thm.cramer.C.1}), applied to
}u=q\right) \\
&  =\left(  a_{1,q},a_{2,q},\ldots,a_{n,q}\right)  ^{T}%
\ \ \ \ \ \ \ \ \ \ \left(  \text{by (\ref{pf.thm.cramer.u-th-col-A}), applied
to }u=q\right)  .
\end{align*}
Thus,%
\[
\left(  c_{1,q},c_{2,q},\ldots,c_{n,q}\right)  =\left(  \underbrace{\left(
c_{1,q},c_{2,q},\ldots,c_{n,q}\right)  ^{T}}_{=\left(  a_{1,q},a_{2,q}%
,\ldots,a_{n,q}\right)  ^{T}}\right)  ^{T}=\left(  \left(  a_{1,q}%
,a_{2,q},\ldots,a_{n,q}\right)  ^{T}\right)  ^{T}=\left(  a_{1,q}%
,a_{2,q},\ldots,a_{n,q}\right)  .
\]
In other words, $c_{p,q}=a_{p,q}$ for every $p\in\left\{  1,2,\ldots
,n\right\}  $. This proves (\ref{pf.thm.cramer.cpq}).}. Now, it is easy to see
that%
\begin{equation}
C_{\sim p,\sim j}=A_{\sim p,\sim j}\ \ \ \ \ \ \ \ \ \ \text{for every }%
p\in\left\{  1,2,\ldots,n\right\}  \label{pf.thm.cramer.CA}%
\end{equation}
\footnote{\textit{Proof of (\ref{pf.thm.cramer.CA}):} Let $p\in\left\{
1,2,\ldots,n\right\}  $.
\par
Let $\left(  u_{1},u_{2},\ldots,u_{n-1}\right)  $ denote the $\left(
n-1\right)  $-tuple $\left(  1,2,\ldots,\widehat{j},\ldots,n\right)  $. Thus,
$\left(  u_{1},u_{2},\ldots,u_{n-1}\right)  =\left(  1,2,\ldots,\widehat{j}%
,\ldots,n\right)  $, so that $\left\{  u_{1},u_{2},\ldots,u_{n-1}\right\}
=\left\{  1,2,\ldots,\widehat{j},\ldots,n\right\}  =\left\{  1,2,\ldots
,n\right\}  \setminus\left\{  j\right\}  $.
\par
Now, let $y\in\left\{  1,2,\ldots,n-1\right\}  $. Then, $u_{y}\in\left\{
u_{1},u_{2},\ldots,u_{n-1}\right\}  =\left\{  1,2,\ldots,n\right\}
\setminus\left\{  j\right\}  $, so that $u_{y}\neq j$. Hence,%
\begin{equation}
c_{p,u_{y}}=a_{p,u_{y}}\ \ \ \ \ \ \ \ \ \ \text{for every }p\in\left\{
1,2,\ldots,n\right\}  \label{pf.thm.cramer.CA.pf.1}%
\end{equation}
(by (\ref{pf.thm.cramer.cpq}), applied to $q=u_{y}$).
\par
Let us now forget that we fixed $y$. We thus have shown that
(\ref{pf.thm.cramer.CA.pf.1}) holds for every $y\in\left\{  1,2,\ldots
,n-1\right\}  $.
\par
Let $\left(  v_{1},v_{2},\ldots,v_{n-1}\right)  $ denote the $\left(
n-1\right)  $-tuple $\left(  1,2,\ldots,\widehat{p},\ldots,n\right)  $. Thus,
$\left(  v_{1},v_{2},\ldots,v_{n-1}\right)  =\left(  1,2,\ldots,\widehat{p}%
,\ldots,n\right)  $.
\par
Now, the definition of $C_{\sim p,\sim j}$ yields%
\begin{align*}
C_{\sim p,\sim j}  &  =\operatorname*{sub}\nolimits_{1,2,\ldots,\widehat{p}%
,\ldots,n}^{1,2,\ldots,\widehat{j},\ldots,n}C=\operatorname*{sub}%
\nolimits_{1,2,\ldots,\widehat{p},\ldots,n}^{u_{1},u_{2},\ldots,u_{n-1}}C\\
&  \ \ \ \ \ \ \ \ \ \ \left(  \text{since }\left(  1,2,\ldots,\widehat{j}%
,\ldots,n\right)  =\left(  u_{1},u_{2},\ldots,u_{n-1}\right)  \right) \\
&  =\operatorname*{sub}\nolimits_{v_{1},v_{2},\ldots,v_{n-1}}^{u_{1}%
,u_{2},\ldots,u_{n-1}}C\ \ \ \ \ \ \ \ \ \ \left(  \text{since }\left(
1,2,\ldots,\widehat{p},\ldots,n\right)  =\left(  v_{1},v_{2},\ldots
,v_{n-1}\right)  \right) \\
&  =\left(  \underbrace{c_{v_{x},u_{y}}}_{\substack{=a_{v_{x},u_{y}%
}\\\text{(by (\ref{pf.thm.cramer.CA.pf.1}),}\\\text{applied to }%
p=v_{x}\text{)}}}\right)  _{1\leq x\leq n-1,\ 1\leq y\leq n-1}\\
&  \ \ \ \ \ \ \ \ \ \ \left(  \text{by the definition of }\operatorname*{sub}%
\nolimits_{v_{1},v_{2},\ldots,v_{n-1}}^{u_{1},u_{2},\ldots,u_{n-1}}C\text{,
since }C=\left(  c_{i,j}\right)  _{1\leq i\leq n,\ 1\leq j\leq n}\right) \\
&  =\left(  a_{v_{x},u_{y}}\right)  _{1\leq x\leq n-1,\ 1\leq y\leq n-1}.
\end{align*}
Compared with%
\begin{align*}
A_{\sim p,\sim j}  &  =\operatorname*{sub}\nolimits_{1,2,\ldots,\widehat{p}%
,\ldots,n}^{1,2,\ldots,\widehat{j},\ldots,n}A=\operatorname*{sub}%
\nolimits_{1,2,\ldots,\widehat{p},\ldots,n}^{u_{1},u_{2},\ldots,u_{n-1}}A\\
&  \ \ \ \ \ \ \ \ \ \ \left(  \text{since }\left(  1,2,\ldots,\widehat{j}%
,\ldots,n\right)  =\left(  u_{1},u_{2},\ldots,u_{n-1}\right)  \right) \\
&  =\operatorname*{sub}\nolimits_{v_{1},v_{2},\ldots,v_{n-1}}^{u_{1}%
,u_{2},\ldots,u_{n-1}}A\ \ \ \ \ \ \ \ \ \ \left(  \text{since }\left(
1,2,\ldots,\widehat{p},\ldots,n\right)  =\left(  v_{1},v_{2},\ldots
,v_{n-1}\right)  \right) \\
&  =\left(  a_{v_{x},u_{y}}\right)  _{1\leq x\leq n-1,\ 1\leq y\leq n-1}\\
&  \ \ \ \ \ \ \ \ \ \ \left(  \text{by the definition of }\operatorname*{sub}%
\nolimits_{v_{1},v_{2},\ldots,v_{n-1}}^{u_{1},u_{2},\ldots,u_{n-1}}A\text{,
since }A=\left(  a_{i,j}\right)  _{1\leq i\leq n,\ 1\leq j\leq n}\right)  ,
\end{align*}
this yields $C_{\sim p,\sim j}=A_{\sim p,\sim j}$. This proves
(\ref{pf.thm.cramer.CA}).}.

Now,%
\begin{align}
\det\left(  \underbrace{A_{j}^{\#}}_{=C}\right)   &  =\det C=\sum_{p=1}%
^{n}\left(  -1\right)  ^{p+j}\underbrace{c_{p,j}}_{\substack{=b_{p}\\\text{(by
(\ref{pf.thm.cramer.cpj}))}}}\det\left(  \underbrace{C_{\sim p,\sim j}%
}_{\substack{=A_{\sim p,\sim j}\\\text{(by (\ref{pf.thm.cramer.CA}))}}}\right)
\nonumber\\
&  \ \ \ \ \ \ \ \ \ \ \left(
\begin{array}
[c]{c}%
\text{by Theorem \ref{thm.laplace.gen} \textbf{(b)}, applied to }C\text{,
}c_{i,j}\text{ and }j\\
\text{instead of }A\text{, }a_{i,j}\text{ and }q
\end{array}
\right) \nonumber\\
&  =\sum_{p=1}^{n}\left(  -1\right)  ^{p+j}b_{p}\det\left(  A_{\sim p,\sim
j}\right)  . \label{pf.thm.cramer.detAshj}%
\end{align}

Let us now forget that we fixed $j$. We thus have proven
(\ref{pf.thm.cramer.detAshj}) for every $j\in\left\{  1,2,\ldots,n\right\}  $.

Now, fix $i\in\left\{  1,2,\ldots,n\right\}  $. Then, for every $p\in\left\{
1,2,\ldots,n\right\}  $ satisfying $p\neq i$, we have%
\begin{equation}
\sum_{q=1}^{n}b_{p}\left(  -1\right)  ^{p+q}a_{i,q}\det\left(  A_{\sim p,\sim
q}\right)  =0 \label{pf.thm.cramer.termkiller1}%
\end{equation}
\footnote{\textit{Proof of (\ref{pf.thm.cramer.termkiller1}):} Let
$p\in\left\{  1,2,\ldots,n\right\}  $ be such that $p\neq i$. Hence,
Proposition \ref{prop.laplace.0} \textbf{(a)} (applied to $r=i$) shows that%
\begin{equation}
0=\sum_{q=1}^{n}\left(  -1\right)  ^{p+q}a_{i,q}\det\left(  A_{\sim p,\sim
q}\right)  . \label{pf.thm.cramer.termkiller1.pf.1}%
\end{equation}
Now,%
\[
\sum_{q=1}^{n}b_{p}\left(  -1\right)  ^{p+q}a_{i,q}\det\left(  A_{\sim p,\sim
q}\right)  =b_{p}\underbrace{\sum_{q=1}^{n}\left(  -1\right)  ^{p+q}%
a_{i,q}\det\left(  A_{\sim p,\sim q}\right)  }_{\substack{=0\\\text{(by
(\ref{pf.thm.cramer.termkiller1.pf.1}))}}}=0.
\]
Thus, (\ref{pf.thm.cramer.termkiller1}) is proven.}. Also, we have%
\begin{equation}
\sum_{q=1}^{n}b_{i}\left(  -1\right)  ^{i+q}a_{i,q}\det\left(  A_{\sim i,\sim
q}\right)  =\det A\cdot b_{i} \label{pf.thm.cramer.termkiller2}%
\end{equation}
\footnote{\textit{Proof of (\ref{pf.thm.cramer.termkiller2}):} Applying
Theorem \ref{thm.laplace.gen} \textbf{(a)} to $p=i$, we obtain%
\begin{equation}
\det A=\sum_{q=1}^{n}\left(  -1\right)  ^{i+q}a_{i,q}\det\left(  A_{\sim
i,\sim q}\right)  . \label{pf.thm.cramer.termkiller2.pf.1}%
\end{equation}
Now,%
\[
\sum_{q=1}^{n}b_{i}\left(  -1\right)  ^{i+q}a_{i,q}\det\left(  A_{\sim i,\sim
q}\right)  =b_{i}\underbrace{\sum_{q=1}^{n}\left(  -1\right)  ^{i+q}%
a_{i,q}\det\left(  A_{\sim i,\sim q}\right)  }_{\substack{=\det A\\\text{(by
(\ref{pf.thm.cramer.termkiller2.pf.1}))}}}=b_{i}\det A=\det A\cdot b_{i}.
\]
This proves (\ref{pf.thm.cramer.termkiller2}).}. Now,%

\begin{align}
&  \sum_{k=1}^{n}a_{i,k}\det\left(  A_{k}^{\#}\right) \nonumber\\
&  =\sum_{j=1}^{n}a_{i,j}\underbrace{\det\left(  A_{j}^{\#}\right)
}_{\substack{=\sum_{p=1}^{n}\left(  -1\right)  ^{p+j}b_{p}\det\left(  A_{\sim
p,\sim j}\right)  \\\text{(by (\ref{pf.thm.cramer.detAshj}))}}%
}\ \ \ \ \ \ \ \ \ \ \left(  \text{here, we renamed the summation index
}k\text{ as }j\right) \nonumber\\
&  =\sum_{j=1}^{n}a_{i,j}\sum_{p=1}^{n}\left(  -1\right)  ^{p+j}b_{p}%
\det\left(  A_{\sim p,\sim j}\right)  =\sum_{q=1}^{n}a_{i,q}\sum_{p=1}%
^{n}\left(  -1\right)  ^{p+q}b_{p}\det\left(  A_{\sim p,\sim q}\right)
\nonumber\\
&  \ \ \ \ \ \ \ \ \ \ \left(  \text{here, we renamed the summation index
}j\text{ as }q\text{ in the first sum}\right) \nonumber\\
&  =\underbrace{\sum_{q=1}^{n}\sum_{p=1}^{n}}_{=\sum_{p=1}^{n}\sum_{q=1}^{n}%
}\underbrace{a_{i,q}\left(  -1\right)  ^{p+q}b_{p}}_{=b_{p}\left(  -1\right)
^{p+q}a_{i,q}}\det\left(  A_{\sim p,\sim q}\right) \nonumber\\
&  =\underbrace{\sum_{p=1}^{n}}_{=\sum_{p\in\left\{  1,2,\ldots,n\right\}  }%
}\sum_{q=1}^{n}b_{p}\left(  -1\right)  ^{p+q}a_{i,q}\det\left(  A_{\sim p,\sim
q}\right) \nonumber\\
&  =\sum_{p\in\left\{  1,2,\ldots,n\right\}  }\sum_{q=1}^{n}b_{p}\left(
-1\right)  ^{p+q}a_{i,q}\det\left(  A_{\sim p,\sim q}\right) \nonumber\\
&  =\sum_{\substack{p\in\left\{  1,2,\ldots,n\right\}  ;\\p\neq i}%
}\underbrace{\sum_{q=1}^{n}b_{p}\left(  -1\right)  ^{p+q}a_{i,q}\det\left(
A_{\sim p,\sim q}\right)  }_{\substack{=0\\\text{(by
(\ref{pf.thm.cramer.termkiller1}))}}}+\underbrace{\sum_{q=1}^{n}b_{i}\left(
-1\right)  ^{i+q}a_{i,q}\det\left(  A_{\sim i,\sim q}\right)  }%
_{\substack{=\det A\cdot b_{i}\\\text{(by (\ref{pf.thm.cramer.termkiller2}))}%
}}\nonumber\\
&  \ \ \ \ \ \ \ \ \ \ \left(  \text{here, we have split off the addend for
}p=i\text{ from the sum}\right) \nonumber\\
&  =\underbrace{\sum_{\substack{p\in\left\{  1,2,\ldots,n\right\}  ;\\p\neq
i}}0}_{=0}+\det A\cdot b_{i}=\det A\cdot b_{i}.
\label{pf.thm.cramer.almostthere}%
\end{align}

Now, let us forget that we fixed $i$. We thus have proven
(\ref{pf.thm.cramer.almostthere}) for every $i\in\left\{  1,2,\ldots
,n\right\}  $.

Now, let $d$ be the vector $\left(  \det\left(  A_{1}^{\#}\right)
,\det\left(  A_{2}^{\#}\right)  ,\ldots,\det\left(  A_{n}^{\#}\right)
\right)  ^{T}$. Thus,%
\begin{align*}
d  &  =\left(  \det\left(  A_{1}^{\#}\right)  ,\det\left(  A_{2}^{\#}\right)
,\ldots,\det\left(  A_{n}^{\#}\right)  \right)  ^{T}=\left(
\begin{array}
[c]{c}%
\det\left(  A_{1}^{\#}\right) \\
\det\left(  A_{2}^{\#}\right) \\
\vdots\\
\det\left(  A_{n}^{\#}\right)
\end{array}
\right) \\
&  =\left(  \det\left(  A_{i}^{\#}\right)  \right)  _{1\leq i\leq n,\ 1\leq
j\leq1}.
\end{align*}

The definition of the product of two matrices shows that%
\begin{align*}
A\cdot d  &  =\left(  \underbrace{\sum_{k=1}^{n}a_{i,k}\det\left(  A_{k}%
^{\#}\right)  }_{\substack{=\det A\cdot b_{i}\\\text{(by
(\ref{pf.thm.cramer.almostthere}))}}}\right)  _{1\leq i\leq n,\ 1\leq j\leq
1}\\
&  \ \ \ \ \ \ \ \ \ \ \left(  \text{since }A=\left(  a_{i,j}\right)  _{1\leq
i\leq n,\ 1\leq j\leq n}\text{ and }d=\left(  \det\left(  A_{i}^{\#}\right)
\right)  _{1\leq i\leq n,\ 1\leq j\leq1}\right) \\
&  =\left(  \det A\cdot b_{i}\right)  _{1\leq i\leq n,\ 1\leq j\leq1}.
\end{align*}
Comparing this with%
\begin{align*}
\det A\cdot\underbrace{b}_{=\left(  b_{1},b_{2},\ldots,b_{n}\right)  ^{T}}  &
=\det A\cdot\left(  b_{1},b_{2},\ldots,b_{n}\right)  ^{T}\\
&  =\left(  \underbrace{\det A\cdot\left(  b_{1},b_{2},\ldots,b_{n}\right)
}_{\substack{=\left(  \det A\cdot b_{1},\det A\cdot b_{2},\ldots,\det A\cdot
b_{n}\right)  \\=\left(  \det A\cdot b_{j}\right)  _{1\leq i\leq1,\ 1\leq
j\leq n}}}\right)  ^{T}\\
&  =\left(  \left(  \det A\cdot b_{j}\right)  _{1\leq i\leq1,\ 1\leq j\leq
n}\right)  ^{T}=\left(  \det A\cdot b_{i}\right)  _{1\leq i\leq n,\ 1\leq
j\leq1}\\
&  \ \ \ \ \ \ \ \ \ \ \left(  \text{by the definition of the transpose of a
matrix}\right)  ,
\end{align*}
we obtain $A\cdot d=\det A\cdot b$. Since $d=\left(  \det\left(  A_{1}%
^{\#}\right)  ,\det\left(  A_{2}^{\#}\right)  ,\ldots,\det\left(  A_{n}%
^{\#}\right)  \right)  ^{T}$, we can rewrite this as $A\cdot\left(
\det\left(  A_{1}^{\#}\right)  ,\det\left(  A_{2}^{\#}\right)  ,\ldots
,\det\left(  A_{n}^{\#}\right)  \right)  ^{T}=\det A\cdot b$. This proves
Theorem \ref{thm.cramer} \textbf{(a)}.

\textbf{(b)} Theorem \ref{thm.matrices.inverses.square} \textbf{(a)} shows
that the matrix $A$ is invertible if and only if the element $\det A$ of
$\mathbb{K}$ is invertible (in $\mathbb{K}$). Hence, the element $\det A$ of
$\mathbb{K}$ is invertible (since the matrix $A$ is invertible). Thus,
$\dfrac{1}{\det A}$ is well-defined. Clearly, $\underbrace{\dfrac{1}{\det
A}\cdot\det A}_{=1}\cdot b=b$, so that%
\begin{align*}
b  &  =\dfrac{1}{\det A}\cdot\underbrace{\det A\cdot b}_{\substack{=A\cdot
\left(  \det\left(  A_{1}^{\#}\right)  ,\det\left(  A_{2}^{\#}\right)
,\ldots,\det\left(  A_{n}^{\#}\right)  \right)  ^{T}\\\text{(by Theorem
\ref{thm.cramer} \textbf{(a)})}}}\\
&  =\dfrac{1}{\det A}\cdot A\cdot\left(  \det\left(  A_{1}^{\#}\right)
,\det\left(  A_{2}^{\#}\right)  ,\ldots,\det\left(  A_{n}^{\#}\right)
\right)  ^{T}\\
&  =A\cdot\underbrace{\left(  \dfrac{1}{\det A}\cdot\left(  \det\left(
A_{1}^{\#}\right)  ,\det\left(  A_{2}^{\#}\right)  ,\ldots,\det\left(
A_{n}^{\#}\right)  \right)  ^{T}\right)  }_{\substack{=\left(  \dfrac{1}{\det
A}\det\left(  A_{1}^{\#}\right)  ,\dfrac{1}{\det A}\det\left(  A_{2}%
^{\#}\right)  ,\ldots,\dfrac{1}{\det A}\det\left(  A_{n}^{\#}\right)  \right)
^{T}\\=\left(  \dfrac{\det\left(  A_{1}^{\#}\right)  }{\det A},\dfrac
{\det\left(  A_{2}^{\#}\right)  }{\det A},\ldots,\dfrac{\det\left(  A_{n}%
^{\#}\right)  }{\det A}\right)  ^{T}}}\\
&  =A\cdot\left(  \dfrac{\det\left(  A_{1}^{\#}\right)  }{\det A},\dfrac
{\det\left(  A_{2}^{\#}\right)  }{\det A},\ldots,\dfrac{\det\left(  A_{n}%
^{\#}\right)  }{\det A}\right)  ^{T}.
\end{align*}
Therefore,
\begin{align*}
&  A^{-1}\underbrace{b}_{=A\cdot\left(  \dfrac{\det\left(  A_{1}^{\#}\right)
}{\det A},\dfrac{\det\left(  A_{2}^{\#}\right)  }{\det A},\ldots,\dfrac
{\det\left(  A_{n}^{\#}\right)  }{\det A}\right)  ^{T}}\\
&  =\underbrace{A^{-1}A}_{=I_{n}}\cdot\left(  \dfrac{\det\left(  A_{1}%
^{\#}\right)  }{\det A},\dfrac{\det\left(  A_{2}^{\#}\right)  }{\det A}%
,\ldots,\dfrac{\det\left(  A_{n}^{\#}\right)  }{\det A}\right)  ^{T}\\
&  =I_{n}\cdot\left(  \dfrac{\det\left(  A_{1}^{\#}\right)  }{\det A}%
,\dfrac{\det\left(  A_{2}^{\#}\right)  }{\det A},\ldots,\dfrac{\det\left(
A_{n}^{\#}\right)  }{\det A}\right)  ^{T}=\left(  \dfrac{\det\left(
A_{1}^{\#}\right)  }{\det A},\dfrac{\det\left(  A_{2}^{\#}\right)  }{\det
A},\ldots,\dfrac{\det\left(  A_{n}^{\#}\right)  }{\det A}\right)  ^{T}.
\end{align*}
This proves Theorem \ref{thm.cramer} \textbf{(b)}.
\end{proof}
\end{verlong}

\subsection{\label{sect.desnanot}The Desnanot-Jacobi identity}

We now move towards more exotic places. In this section\footnote{which,
unfortunately, has become more technical and tedious than it promised to be --
for which I apologize}, we shall prove the \textit{Desnanot-Jacobi identity},
also known as \textit{Lewis Carroll identity}\footnote{See \cite[\S 3.5]%
{Bresso99} for the history of this identity (as well as for a proof different
from ours, and for an application). In a nutshell: Desnanot discovered it in
1819; Jacobi proved it in 1833 (and again in 1841); in 1866, Charles Lutwidge
Dodgson (better known as Lewis Carroll, although his mathematical works were
not printed under this pen name) popularized it by publishing an algorithm for
evaluating determinants that made heavy use of this identity.}. We will need
some notations to state this identity in the generality I want; but first I
shall state the two best known particular cases (from which the general
version can actually be easily derived, although this is not the way I will
take):\footnote{We shall use the notations of Definition
\ref{def.submatrix.minor} throughout this section.}

\begin{proposition}
\label{prop.desnanot.1n}Let $n\in\mathbb{N}$ be such that $n\geq2$. Let
$A=\left(  a_{i,j}\right)  _{1\leq i\leq n,\ 1\leq j\leq n}$ be an $n\times
n$-matrix. Let $A^{\prime}$ be the $\left(  n-2\right)  \times\left(
n-2\right)  $-matrix $\left(  a_{i+1,j+1}\right)  _{1\leq i\leq n-2,\ 1\leq
j\leq n-2}$. (In other words, $A^{\prime}$ is what remains of the matrix $A$
when we remove the first row, the last row, the first column and the last
column.) Then,%
\begin{align*}
&  \det A\cdot\det\left(  A^{\prime}\right) \\
&  =\det\left(  A_{\sim1,\sim1}\right)  \cdot\det\left(  A_{\sim n,\sim
n}\right)  -\det\left(  A_{\sim1,\sim n}\right)  \cdot\det\left(  A_{\sim
n,\sim1}\right)  .
\end{align*}

\end{proposition}

\begin{proposition}
\label{prop.desnanot.12}Let $n\in\mathbb{N}$ be such that $n\geq2$. Let
$A=\left(  a_{i,j}\right)  _{1\leq i\leq n,\ 1\leq j\leq n}$ be an $n\times
n$-matrix. Let $\widetilde{A}$ be the $\left(  n-2\right)  \times\left(
n-2\right)  $-matrix $\left(  a_{i+2,j+2}\right)  _{1\leq i\leq n-2,\ 1\leq
j\leq n-2}$. (In other words, $\widetilde{A}$ is what remains of the matrix
$A$ when we remove the first two rows and the first two columns.) Then,%
\begin{align*}
&  \det A\cdot\det\widetilde{A}\\
&  =\det\left(  A_{\sim1,\sim1}\right)  \cdot\det\left(  A_{\sim2,\sim
2}\right)  -\det\left(  A_{\sim1,\sim2}\right)  \cdot\det\left(  A_{\sim
2,\sim1}\right)  .
\end{align*}

\end{proposition}

\begin{example}
For this example, set $n=4$ and $A=\left(
\begin{array}
[c]{cccc}%
a_{1} & b_{1} & c_{1} & d_{1}\\
a_{2} & b_{2} & c_{2} & d_{2}\\
a_{3} & b_{3} & c_{3} & d_{3}\\
a_{4} & b_{4} & c_{4} & d_{4}%
\end{array}
\right)  $. Then, Proposition \ref{prop.desnanot.1n} says that%
\begin{align*}
&  \det\left(
\begin{array}
[c]{cccc}%
a_{1} & b_{1} & c_{1} & d_{1}\\
a_{2} & b_{2} & c_{2} & d_{2}\\
a_{3} & b_{3} & c_{3} & d_{3}\\
a_{4} & b_{4} & c_{4} & d_{4}%
\end{array}
\right)  \cdot\det\left(
\begin{array}
[c]{cc}%
b_{2} & c_{2}\\
b_{3} & c_{3}%
\end{array}
\right) \\
&  =\det\left(
\begin{array}
[c]{ccc}%
b_{2} & c_{2} & d_{2}\\
b_{3} & c_{3} & d_{3}\\
b_{4} & c_{4} & d_{4}%
\end{array}
\right)  \cdot\det\left(
\begin{array}
[c]{ccc}%
a_{1} & b_{1} & c_{1}\\
a_{2} & b_{2} & c_{2}\\
a_{3} & b_{3} & c_{3}%
\end{array}
\right) \\
&  \ \ \ \ \ \ \ \ \ \ -\det\left(
\begin{array}
[c]{ccc}%
a_{2} & b_{2} & c_{2}\\
a_{3} & b_{3} & c_{3}\\
a_{4} & b_{4} & c_{4}%
\end{array}
\right)  \cdot\det\left(
\begin{array}
[c]{ccc}%
b_{1} & c_{1} & d_{1}\\
b_{2} & c_{2} & d_{2}\\
b_{3} & c_{3} & d_{3}%
\end{array}
\right)  .
\end{align*}
Meanwhile, Proposition \ref{prop.desnanot.12} says that%
\begin{align*}
&  \det\left(
\begin{array}
[c]{cccc}%
a_{1} & b_{1} & c_{1} & d_{1}\\
a_{2} & b_{2} & c_{2} & d_{2}\\
a_{3} & b_{3} & c_{3} & d_{3}\\
a_{4} & b_{4} & c_{4} & d_{4}%
\end{array}
\right)  \cdot\det\left(
\begin{array}
[c]{cc}%
c_{3} & d_{3}\\
c_{4} & d_{4}%
\end{array}
\right) \\
&  =\det\left(
\begin{array}
[c]{ccc}%
b_{2} & c_{2} & d_{2}\\
b_{3} & c_{3} & d_{3}\\
b_{4} & c_{4} & d_{4}%
\end{array}
\right)  \cdot\det\left(
\begin{array}
[c]{ccc}%
a_{1} & c_{1} & d_{1}\\
a_{3} & c_{3} & d_{3}\\
a_{4} & c_{4} & d_{4}%
\end{array}
\right) \\
&  \ \ \ \ \ \ \ \ \ \ -\det\left(
\begin{array}
[c]{ccc}%
a_{2} & c_{2} & d_{2}\\
a_{3} & c_{3} & d_{3}\\
a_{4} & c_{4} & d_{4}%
\end{array}
\right)  \cdot\det\left(
\begin{array}
[c]{ccc}%
b_{1} & c_{1} & d_{1}\\
b_{3} & c_{3} & d_{3}\\
b_{4} & c_{4} & d_{4}%
\end{array}
\right)  .
\end{align*}

\end{example}

Proposition \ref{prop.desnanot.1n} occurs (for instance) in
\cite[\textit{(Alice)}]{zeilberger-twotime}, in \cite[Theorem 3.12]{Bresso99},
in \cite[Example 3.3]{Gill15} and in \cite[Proposition 10]{Krattenthaler}
(without a proof, but with a brief list of applications). Proposition
\ref{prop.desnanot.12} occurs (among other places) in \cite[(1)]{BerBru08}
(with a generalization). The reader can easily see that Proposition
\ref{prop.desnanot.12} is equivalent to Proposition \ref{prop.desnanot.1n}%
\footnote{Indeed, one is easily obtained from the other by swapping the $2$-nd
and the $n$-th rows of the matrix $A$ and swapping the $2$-nd and the $n$-th
columns of the matrix $A$, applying parts \textbf{(a)} and \textbf{(b)} of
Exercise \ref{exe.ps4.6} and checking that all signs cancel.}; we shall prove
a generalization of both.

Let me now introduce some notations:\footnote{Recall that we are using the
notations of Definition \ref{def.rowscols}, of Definition \ref{def.hat-omit},
and of Definition \ref{def.submatrix.minor}.}

\begin{definition}
\label{def.unrows2}Let $n\in\mathbb{N}$. Let $r$ and $s$ be two elements of
$\left\{  1,2,\ldots,n\right\}  $ such that $r<s$. Then, $\left(
1,2,\ldots,\widehat{r},\ldots,\widehat{s},\ldots,n\right)  $ will denote the
$\left(  n-2\right)  $-tuple
\[
\left(  \underbrace{1,2,\ldots,r-1}_{\substack{\text{all integers}\\\text{from
}1\text{ to }r-1}},\underbrace{r+1,r+2,\ldots,s-1}_{\substack{\text{all
integers}\\\text{from }r+1\text{ to }s-1}},\underbrace{s+1,s+2,\ldots
,n}_{\substack{\text{all integers}\\\text{from }s+1\text{ to }n}}\right)  .
\]
In other words, $\left(  1,2,\ldots,\widehat{r},\ldots,\widehat{s}%
,\ldots,n\right)  $ will denote the result of removing the entries $r$ and $s$
from the $n$-tuple $\left(  1,2,\ldots,n\right)  $.
\end{definition}

We can now state a more general version of the Desnanot-Jacobi identity:

\begin{theorem}
\label{thm.desnanot}Let $n\in\mathbb{N}$ be such that $n\geq2$. Let $p$, $q$,
$u$ and $v$ be four elements of $\left\{  1,2,\ldots,n\right\}  $ such that
$p<q$ and $u<v$. Let $A$ be an $n\times n$-matrix. Then,%
\begin{align*}
&  \det A\cdot\det\left(  \operatorname*{sub}\nolimits_{1,2,\ldots
,\widehat{p},\ldots,\widehat{q},\ldots,n}^{1,2,\ldots,\widehat{u}%
,\ldots,\widehat{v},\ldots,n}A\right) \\
&  =\det\left(  A_{\sim p,\sim u}\right)  \cdot\det\left(  A_{\sim q,\sim
v}\right)  -\det\left(  A_{\sim p,\sim v}\right)  \cdot\det\left(  A_{\sim
q,\sim u}\right)  .
\end{align*}

\end{theorem}

\begin{example}
If we set $n=3$, $p=1$, $q=2$, $u=1$, $v=3$ and $A=\left(
\begin{array}
[c]{ccc}%
a & a^{\prime} & a^{\prime\prime}\\
b & b^{\prime} & b^{\prime\prime}\\
c & c^{\prime} & c^{\prime\prime}%
\end{array}
\right)  $, then Theorem \ref{thm.desnanot} says that%
\begin{align*}
&  \det\left(
\begin{array}
[c]{ccc}%
a & a^{\prime} & a^{\prime\prime}\\
b & b^{\prime} & b^{\prime\prime}\\
c & c^{\prime} & c^{\prime\prime}%
\end{array}
\right)  \det\left(
\begin{array}
[c]{c}%
c^{\prime}%
\end{array}
\right) \\
&  =\det\left(
\begin{array}
[c]{cc}%
b^{\prime} & b^{\prime\prime}\\
c^{\prime} & c^{\prime\prime}%
\end{array}
\right)  \cdot\det\left(
\begin{array}
[c]{cc}%
a & a^{\prime}\\
c & c^{\prime}%
\end{array}
\right)  -\det\left(
\begin{array}
[c]{cc}%
b & b^{\prime}\\
c & c^{\prime}%
\end{array}
\right)  \cdot\det\left(
\begin{array}
[c]{cc}%
a^{\prime} & a^{\prime\prime}\\
c^{\prime} & c^{\prime\prime}%
\end{array}
\right)  .
\end{align*}

\end{example}

Before we prove this theorem, let me introduce some more notations:

\begin{definition}
\label{def.unrows}Let $n\in\mathbb{N}$ and $m\in\mathbb{N}$. Let
$A\in\mathbb{K}^{n\times m}$ be an $n\times m$-matrix.

\textbf{(a)} For every $u\in\left\{  1,2,\ldots,n\right\}  $, we let
$A_{u,\bullet}$ be the $u$-th row of the matrix $A$. This $A_{u,\bullet}$ is
thus a row vector with $m$ entries, i.e., a $1\times m$-matrix.

\textbf{(b)} For every $v\in\left\{  1,2,\ldots,m\right\}  $, we let
$A_{\bullet,v}$ be the $v$-th column of the matrix $A$. This $A_{\bullet,v}$
is thus a column vector with $n$ entries, i.e., an $n\times1$-matrix.

\textbf{(c)} For every $u\in\left\{  1,2,\ldots,n\right\}  $, we set $A_{\sim
u,\bullet}=\operatorname*{rows}\nolimits_{1,2,\ldots,\widehat{u},\ldots,n}A$.
This $A_{\sim u,\bullet}$ is thus an $\left(  n-1\right)  \times m$-matrix.
(In more intuitive terms, the definition of $A_{\sim u,\bullet}$ rewrites as
follows: $A_{\sim u,\bullet}$ is the matrix obtained from the matrix $A$ by
removing the $u$-th row.)

\textbf{(d)} For every $v\in\left\{  1,2,\ldots,m\right\}  $, we set
$A_{\bullet,\sim v}=\operatorname*{cols}\nolimits_{1,2,\ldots,\widehat{v}%
,\ldots,m}A$. This $A_{\bullet,\sim v}$ is thus an $n\times\left(  m-1\right)
$-matrix. (In more intuitive terms, the definition of $A_{\bullet,\sim v}$
rewrites as follows: $A_{\bullet,\sim v}$ is the matrix obtained from the
matrix $A$ by removing the $v$-th column.)
\end{definition}

\begin{example}
\label{exam.unrows}If $n=3$, $m=4$ and $A=\left(
\begin{array}
[c]{cccc}%
a & b & c & d\\
a^{\prime} & b^{\prime} & c^{\prime} & d^{\prime}\\
a^{\prime\prime} & b^{\prime\prime} & c^{\prime\prime} & d^{\prime\prime}%
\end{array}
\right)  $, then%
\begin{align*}
A_{2,\bullet}  &  =\left(
\begin{array}
[c]{cccc}%
a^{\prime} & b^{\prime} & c^{\prime} & d^{\prime}%
\end{array}
\right)  ,\ \ \ \ \ \ \ \ \ \ A_{\bullet,2}=\left(
\begin{array}
[c]{c}%
b\\
b^{\prime}\\
b^{\prime\prime}%
\end{array}
\right)  ,\\
A_{\sim2,\bullet}  &  =\left(
\begin{array}
[c]{cccc}%
a & b & c & d\\
a^{\prime\prime} & b^{\prime\prime} & c^{\prime\prime} & d^{\prime\prime}%
\end{array}
\right)  ,\ \ \ \ \ \ \ \ \ \ A_{\bullet,\sim2}=\left(
\begin{array}
[c]{ccc}%
a & c & d\\
a^{\prime} & c^{\prime} & d^{\prime}\\
a^{\prime\prime} & c^{\prime\prime} & d^{\prime\prime}%
\end{array}
\right)  .
\end{align*}

\end{example}

Here are some simple properties of the notations introduced in Definition
\ref{def.unrows}:

\begin{proposition}
\label{prop.unrows.basics}Let $n\in\mathbb{N}$ and $m\in\mathbb{N}$. Let
$A\in\mathbb{K}^{n\times m}$ be an $n\times m$-matrix.

\textbf{(a)} For every $u\in\left\{  1,2,\ldots,n\right\}  $, we have%
\[
A_{u,\bullet}=\left(  \text{the }u\text{-th row of the matrix }A\right)
=\operatorname*{rows}\nolimits_{u}A.
\]
(Here, the notation $\operatorname*{rows}\nolimits_{u}A$ is a particular case
of Definition \ref{def.rowscols} \textbf{(a)}.)

\textbf{(b)} For every $v\in\left\{  1,2,\ldots,m\right\}  $, we have%
\[
A_{\bullet,v}=\left(  \text{the }v\text{-th column of the matrix }A\right)
=\operatorname*{cols}\nolimits_{v}A.
\]

\textbf{(c)} For every $u\in\left\{  1,2,\ldots,n\right\}  $ and $v\in\left\{
1,2,\ldots,m\right\}  $, we have%
\[
\left(  A_{\bullet,\sim v}\right)  _{\sim u,\bullet}=\left(  A_{\sim
u,\bullet}\right)  _{\bullet,\sim v}=A_{\sim u,\sim v}.
\]

\textbf{(d)} For every $v\in\left\{  1,2,\ldots,m\right\}  $ and $w\in\left\{
1,2,\ldots,v-1\right\}  $, we have $\left(  A_{\bullet,\sim v}\right)
_{\bullet,w}=A_{\bullet,w}$.

\textbf{(e)} For every $v\in\left\{  1,2,\ldots,m\right\}  $ and $w\in\left\{
v,v+1,\ldots,m-1\right\}  $, we have $\left(  A_{\bullet,\sim v}\right)
_{\bullet,w}=A_{\bullet,w+1}$.

\textbf{(f)} For every $u\in\left\{  1,2,\ldots,n\right\}  $ and $w\in\left\{
1,2,\ldots,u-1\right\}  $, we have $\left(  A_{\sim u,\bullet}\right)
_{w,\bullet}=A_{w,\bullet}$.

\textbf{(g)} For every $u\in\left\{  1,2,\ldots,n\right\}  $ and $w\in\left\{
u,u+1,\ldots,n-1\right\}  $, we have $\left(  A_{\sim u,\bullet}\right)
_{w,\bullet}=A_{w+1,\bullet}$.

\textbf{(h)} For every $v\in\left\{  1,2,\ldots,m\right\}  $ and $w\in\left\{
1,2,\ldots,v-1\right\}  $, we have%
\[
\left(  A_{\bullet,\sim v}\right)  _{\bullet,\sim w}=\operatorname*{cols}%
\nolimits_{1,2,\ldots,\widehat{w},\ldots,\widehat{v},\ldots,m}A.
\]

\textbf{(i)} For every $v\in\left\{  1,2,\ldots,m\right\}  $ and $w\in\left\{
v,v+1,\ldots,m-1\right\}  $, we have%
\[
\left(  A_{\bullet,\sim v}\right)  _{\bullet,\sim w}=\operatorname*{cols}%
\nolimits_{1,2,\ldots,\widehat{v},\ldots,\widehat{w+1},\ldots,m}A.
\]

\textbf{(j)} For every $u\in\left\{  1,2,\ldots,n\right\}  $ and $w\in\left\{
1,2,\ldots,u-1\right\}  $, we have%
\[
\left(  A_{\sim u,\bullet}\right)  _{\sim w,\bullet}=\operatorname*{rows}%
\nolimits_{1,2,\ldots,\widehat{w},\ldots,\widehat{u},\ldots,n}A.
\]

\textbf{(k)} For every $u\in\left\{  1,2,\ldots,n\right\}  $ and $w\in\left\{
u,u+1,\ldots,n-1\right\}  $, we have%
\[
\left(  A_{\sim u,\bullet}\right)  _{\sim w,\bullet}=\operatorname*{rows}%
\nolimits_{1,2,\ldots,\widehat{u},\ldots,\widehat{w+1},\ldots,n}A.
\]

\textbf{(l)} For every $v\in\left\{  1,2,\ldots,n\right\}  $, $u\in\left\{
1,2,\ldots,n\right\}  $ and $q\in\left\{  1,2,\ldots,m\right\}  $ satisfying
$u<v$, we have%
\[
\left(  A_{\sim v,\bullet}\right)  _{\sim u,\sim q}=\operatorname*{rows}%
\nolimits_{1,2,\ldots,\widehat{u},\ldots,\widehat{v},\ldots,n}\left(
A_{\bullet,\sim q}\right)  .
\]

\end{proposition}

\begin{proposition}
\label{prop.unrows.basics-I}Let $n\in\mathbb{N}$. Let $u$ and $v$ be two
elements of $\left\{  1,2,\ldots,n\right\}  $ such that $u<v$. Then, $\left(
\left(  I_{n}\right)  _{\bullet,u}\right)  _{\sim v,\bullet}=\left(
I_{n-1}\right)  _{\bullet,u}$. (Recall that $I_{m}$ denotes the $m\times m$
identity matrix for each $m\in\mathbb{N}$.)
\end{proposition}

\begin{exercise}
\label{exe.unrows.basics}\textbf{(a)} Prove Proposition
\ref{prop.unrows.basics} and Proposition \ref{prop.unrows.basics-I}.

\textbf{(b)} Derive Proposition \ref{prop.desnanot.12} and Proposition
\ref{prop.desnanot.1n} from Theorem \ref{thm.desnanot}.
\end{exercise}

Here comes another piece of notation:

\begin{definition}
\label{def.addcol}Let $n\in\mathbb{N}$ and $m\in\mathbb{N}$. Let
$A\in\mathbb{K}^{n\times m}$ be an $n\times m$-matrix. Let $v\in
\mathbb{K}^{n\times1}$ be a column vector with $n$ entries. Then, $\left(
A\mid v\right)  $ will denote the $n\times\left(  m+1\right)  $-matrix whose
$m+1$ columns are $A_{\bullet,1},A_{\bullet,2},\ldots,A_{\bullet,m},v$ (from
left to right). (Informally speaking, $\left(  A\mid v\right)  $ is the matrix
obtained when the column vector $v$ is \textquotedblleft
attached\textquotedblright\ to $A$ at the right edge.)
\end{definition}

\begin{example}
\label{exam.addcol}We have $\left(  \left(
\begin{array}
[c]{cc}%
a & b\\
c & d
\end{array}
\right)  \mid\left(
\begin{array}
[c]{c}%
p\\
q
\end{array}
\right)  \right)  =\left(
\begin{array}
[c]{ccc}%
a & b & p\\
c & d & q
\end{array}
\right)  $.
\end{example}

The following properties of the notation introduced in Definition
\ref{def.addcol} are not too hard to see (see the solution of Exercise
\ref{exe.prop.addcol.props} for their proofs), and will be used below:

\begin{proposition}
\label{prop.addcol.props1}Let $n\in\mathbb{N}$ and $m\in\mathbb{N}$. Let
$A\in\mathbb{K}^{n\times m}$ be an $n\times m$-matrix. Let $v\in
\mathbb{K}^{n\times1}$ be a column vector with $n$ entries.

\textbf{(a)} Every $q\in\left\{  1,2,\ldots,m\right\}  $ satisfies $\left(
A\mid v\right)  _{\bullet,q}=A_{\bullet,q}$.

\textbf{(b)} We have $\left(  A\mid v\right)  _{\bullet,m+1}=v$.

\textbf{(c)} Every $q\in\left\{  1,2,\ldots,m\right\}  $ satisfies $\left(
A\mid v\right)  _{\bullet,\sim q}=\left(  A_{\bullet,\sim q}\mid v\right)  $.

\textbf{(d)} We have $\left(  A\mid v\right)  _{\bullet,\sim\left(
m+1\right)  }=A$.

\textbf{(e)} We have $\left(  A\mid v\right)  _{\sim p,\bullet}=\left(
A_{\sim p,\bullet}\mid v_{\sim p,\bullet}\right)  $ for every $p\in\left\{
1,2,\ldots,n\right\}  $.

\textbf{(f)} We have $\left(  A\mid v\right)  _{\sim p,\sim\left(  m+1\right)
}=A_{\sim p,\bullet}$ for every $p\in\left\{  1,2,\ldots,n\right\}  $.
\end{proposition}

\begin{proposition}
\label{prop.addcol.props2}Let $n$ be a positive integer. Let $A\in
\mathbb{K}^{n\times\left(  n-1\right)  }$.

\textbf{(a)} For every $v=\left(  v_{1},v_{2},\ldots,v_{n}\right)  ^{T}%
\in\mathbb{K}^{n\times1}$, we have%
\[
\det\left(  A\mid v\right)  =\sum_{i=1}^{n}\left(  -1\right)  ^{n+i}v_{i}%
\det\left(  A_{\sim i,\bullet}\right)  .
\]

\textbf{(b)} For every $p\in\left\{  1,2,\ldots,n\right\}  $, we have
\[
\det\left(  A\mid\left(  I_{n}\right)  _{\bullet,p}\right)  =\left(
-1\right)  ^{n+p}\det\left(  A_{\sim p,\bullet}\right)  .
\]
(Notice that $\left(  I_{n}\right)  _{\bullet,p}$ is the $p$-th column of the
$n\times n$ identity matrix, i.e., the column vector $\left(
\underbrace{0,0,\ldots,0}_{p-1\text{ zeroes}},1,\underbrace{0,0,\ldots
,0}_{n-p\text{ zeroes}}\right)  ^{T}$.)
\end{proposition}

\begin{proposition}
\label{prop.addcol.props3}Let $n\in\mathbb{N}$. Let $A\in\mathbb{K}^{n\times
n}$.

\textbf{(a)} If $n>0$, then $\left(  A_{\bullet,\sim n}\mid A_{\bullet
,n}\right)  =A$.

\textbf{(b)} For every $q\in\left\{  1,2,\ldots,n\right\}  $, we have
$\det\left(  A_{\bullet,\sim q}\mid A_{\bullet,q}\right)  =\left(  -1\right)
^{n+q}\det A$.

\textbf{(c)} If $r$ and $q$ are two elements of $\left\{  1,2,\ldots
,n\right\}  $ satisfying $r\neq q$, then $\det\left(  A_{\bullet,\sim q}\mid
A_{\bullet,r}\right)  =0$.

\textbf{(d)} For every $p\in\left\{  1,2,\ldots,n\right\}  $ and $q\in\left\{
1,2,\ldots,n\right\}  $, we have $\det\left(  A_{\bullet,\sim q}\mid\left(
I_{n}\right)  _{\bullet,p}\right)  =\left(  -1\right)  ^{n+p}\det\left(
A_{\sim p,\sim q}\right)  $.

\textbf{(e)} If $u$ and $v$ are two elements of $\left\{  1,2,\ldots
,n\right\}  $ satisfying $u<v$, and if $r$ is an element of $\left\{
1,2,\ldots,n-1\right\}  $ satisfying $r\neq u$, then $\det\left(
A_{\bullet,\sim u}\mid\left(  A_{\bullet,\sim v}\right)  _{\bullet,r}\right)
=0$.

\textbf{(f)} If $u$ and $v$ are two elements of $\left\{  1,2,\ldots
,n\right\}  $ satisfying $u<v$, then $\left(  -1\right)  ^{u}\det\left(
A_{\bullet,\sim u}\mid\left(  A_{\bullet,\sim v}\right)  _{\bullet,u}\right)
=\left(  -1\right)  ^{n}\det A$.
\end{proposition}

\begin{exercise}
\label{exe.prop.addcol.props}Prove Proposition \ref{prop.addcol.props1},
Proposition \ref{prop.addcol.props2} and Proposition \ref{prop.addcol.props3}.
\end{exercise}

Now, we can state a slightly more interesting identity:

\begin{proposition}
\label{prop.desnanot.AC}Let $n$ be a positive integer. Let $A\in
\mathbb{K}^{n\times\left(  n-1\right)  }$ and $C\in\mathbb{K}^{n\times n}$.
Let $v\in\left\{  1,2,\ldots,n\right\}  $. Then,%
\[
\det\left(  A_{\sim v,\bullet}\right)  \det C=\sum_{q=1}^{n}\left(  -1\right)
^{n+q}\det\left(  A\mid C_{\bullet,q}\right)  \det\left(  C_{\sim v,\sim
q}\right)  .
\]

\end{proposition}

\begin{example}
\label{exam.desnanot.AC}If we set $n=3$, $A=\left(
\begin{array}
[c]{cc}%
a & a^{\prime}\\
b & b^{\prime}\\
c & c^{\prime}%
\end{array}
\right)  $, $C=\left(
\begin{array}
[c]{ccc}%
x & x^{\prime} & x^{\prime\prime}\\
y & y^{\prime} & y^{\prime\prime}\\
z & z^{\prime} & z^{\prime\prime}%
\end{array}
\right)  $ and $v=2$, then Proposition \ref{prop.desnanot.AC} states that%
\begin{align*}
&  \det\left(
\begin{array}
[c]{cc}%
a & a^{\prime}\\
c & c^{\prime}%
\end{array}
\right)  \det\left(
\begin{array}
[c]{ccc}%
x & x^{\prime} & x^{\prime\prime}\\
y & y^{\prime} & y^{\prime\prime}\\
z & z^{\prime} & z^{\prime\prime}%
\end{array}
\right) \\
&  =\det\left(
\begin{array}
[c]{ccc}%
a & a^{\prime} & x\\
b & b^{\prime} & y\\
c & c^{\prime} & z
\end{array}
\right)  \det\left(
\begin{array}
[c]{cc}%
x^{\prime} & x^{\prime\prime}\\
z^{\prime} & z^{\prime\prime}%
\end{array}
\right)  -\det\left(
\begin{array}
[c]{ccc}%
a & a^{\prime} & x^{\prime}\\
b & b^{\prime} & y^{\prime}\\
c & c^{\prime} & z^{\prime}%
\end{array}
\right)  \det\left(
\begin{array}
[c]{cc}%
x & x^{\prime\prime}\\
z & z^{\prime\prime}%
\end{array}
\right) \\
&  \ \ \ \ \ \ \ \ \ \ +\det\left(
\begin{array}
[c]{ccc}%
a & a^{\prime} & x^{\prime\prime}\\
b & b^{\prime} & y^{\prime\prime}\\
c & c^{\prime} & z^{\prime\prime}%
\end{array}
\right)  \det\left(
\begin{array}
[c]{cc}%
x & x^{\prime}\\
z & z^{\prime}%
\end{array}
\right)  .
\end{align*}

\end{example}

\begin{proof}
[Proof of Proposition \ref{prop.desnanot.AC}.]Write the $n\times n$-matrix $C$
in the form $C=\left(  c_{i,j}\right)  _{1\leq i\leq n,\ 1\leq j\leq n}$.

Fix $q\in\left\{  1,2,\ldots,n\right\}  $. Then, $C_{\bullet,q}$ is the $q$-th
column of the matrix $C$ (by the definition of $C_{\bullet,q}$). Thus,%
\begin{align*}
C_{\bullet,q}  &  =\left(  \text{the }q\text{-th column of the matrix
}C\right) \\
&  =\left(
\begin{array}
[c]{c}%
c_{1,q}\\
c_{2,q}\\
\vdots\\
c_{n,q}%
\end{array}
\right)  \ \ \ \ \ \ \ \ \ \ \left(  \text{since }C=\left(  c_{i,j}\right)
_{1\leq i\leq n,\ 1\leq j\leq n}\right) \\
&  =\left(  c_{1,q},c_{2,q},\ldots,c_{n,q}\right)  ^{T}.
\end{align*}
Now, Proposition \ref{prop.addcol.props2} \textbf{(a)} (applied to
$C_{\bullet,q}$ and $c_{i,q}$ instead of $v$ and $v_{i}$) yields
\begin{align}
\det\left(  A\mid C_{\bullet,q}\right)   &  =\sum_{i=1}^{n}\left(  -1\right)
^{n+i}c_{i,q}\det\left(  A_{\sim i,\bullet}\right) \nonumber\\
&  =\sum_{p=1}^{n}\left(  -1\right)  ^{n+p}c_{p,q}\det\left(  A_{\sim
p,\bullet}\right)  \label{pf.prop.desnanot.AC.2}%
\end{align}
(here, we have renamed the summation index $i$ as $p$).

Now, forget that we fixed $q$. We thus have proven
(\ref{pf.prop.desnanot.AC.2}) for each $q\in\left\{  1,2,\ldots,n\right\}  $.

We have%
\begin{equation}
\sum_{q=1}^{n}\left(  -1\right)  ^{v+q}c_{v,q}\det\left(  C_{\sim v,\sim
q}\right)  =\det C \label{pf.prop.desnanot.AC.4}%
\end{equation}
\footnote{\textit{Proof of (\ref{pf.prop.desnanot.AC.4}):} Theorem
\ref{thm.laplace.gen} (applied to $C$, $c_{i,j}$ and $v$ instead of $A$,
$a_{i,j}$ and $p$) yields%
\begin{equation}
\det C=\sum_{q=1}^{n}\left(  -1\right)  ^{v+q}c_{v,q}\det\left(  C_{\sim
v,\sim q}\right)  . \label{pf.prop.desnanot.AC.4.pf.1}%
\end{equation}
This proves (\ref{pf.prop.desnanot.AC.4}).}.

On the other hand, every $p\in\left\{  1,2,\ldots,n\right\}  $ satisfying
$p\neq v$ satisfies%
\begin{equation}
\sum_{q=1}^{n}\left(  -1\right)  ^{p+q}c_{p,q}\det\left(  C_{\sim v,\sim
q}\right)  =0 \label{pf.prop.desnanot.AC.3}%
\end{equation}
\footnote{\textit{Proof of (\ref{pf.prop.desnanot.AC.3}):} Let $p\in\left\{
1,2,\ldots,n\right\}  $ be such that $p\neq v$. Thus, $v\neq p$. Hence,
Proposition \ref{prop.laplace.0} \textbf{(a)} (applied to $C$, $c_{i,j}$, $p$
and $v$ instead of $A$, $a_{i,j}$, $r$ and $p$) yields%
\begin{equation}
0=\sum_{q=1}^{n}\left(  -1\right)  ^{v+q}c_{p,q}\det\left(  C_{\sim v,\sim
q}\right)  . \label{pf.prop.desnanot.AC.3.pf.1}%
\end{equation}
Now, every $q\in\left\{  1,2,\ldots,n\right\}  $ satisfies
\begin{align*}
\left(  -1\right)  ^{p+q}  &  =\left(  -1\right)  ^{\left(  p-v\right)
+\left(  v+q\right)  }\ \ \ \ \ \ \ \ \ \ \left(  \text{since }p+q=\left(
p-v\right)  +\left(  v+q\right)  \right) \\
&  =\left(  -1\right)  ^{p-v}\left(  -1\right)  ^{v+q}.
\end{align*}
Hence,%
\begin{align*}
&  \sum_{q=1}^{n}\underbrace{\left(  -1\right)  ^{p+q}}_{=\left(  -1\right)
^{p-v}\left(  -1\right)  ^{v+q}}c_{p,q}\det\left(  C_{\sim v,\sim q}\right) \\
&  =\sum_{q=1}^{n}\left(  -1\right)  ^{p-v}\left(  -1\right)  ^{v+q}%
c_{p,q}\det\left(  C_{\sim v,\sim q}\right)  =\left(  -1\right)
^{p-v}\underbrace{\sum_{q=1}^{n}\left(  -1\right)  ^{v+q}c_{p,q}\det\left(
C_{\sim v,\sim q}\right)  }_{\substack{=0\\\text{(by
(\ref{pf.prop.desnanot.AC.3.pf.1}))}}}=0.
\end{align*}
This proves (\ref{pf.prop.desnanot.AC.3}).}.

Now,%
\begin{align*}
&  \sum_{q=1}^{n}\left(  -1\right)  ^{n+q}\underbrace{\det\left(  A\mid
C_{\bullet,q}\right)  }_{\substack{=\sum_{p=1}^{n}\left(  -1\right)
^{n+p}c_{p,q}\det\left(  A_{\sim p,\bullet}\right)  \\\text{(by
(\ref{pf.prop.desnanot.AC.2}))}}}\det\left(  C_{\sim v,\sim q}\right) \\
&  =\sum_{q=1}^{n}\left(  -1\right)  ^{n+q}\left(  \sum_{p=1}^{n}\left(
-1\right)  ^{n+p}c_{p,q}\det\left(  A_{\sim p,\bullet}\right)  \right)
\det\left(  C_{\sim v,\sim q}\right) \\
&  =\underbrace{\sum_{q=1}^{n}\sum_{p=1}^{n}}_{=\sum_{p=1}^{n}\sum_{q=1}^{n}%
}\underbrace{\left(  -1\right)  ^{n+q}\left(  -1\right)  ^{n+p}}%
_{\substack{=\left(  -1\right)  ^{\left(  n+q\right)  +\left(  n+p\right)
}=\left(  -1\right)  ^{p+q}\\\text{(since }\left(  n+q\right)  +\left(
n+p\right)  =2n+p+q\equiv p+q\operatorname{mod}2\text{)}}}c_{p,q}\det\left(
A_{\sim p,\bullet}\right)  \det\left(  C_{\sim v,\sim q}\right) \\
&  =\underbrace{\sum_{p=1}^{n}}_{=\sum_{p\in\left\{  1,2,\ldots,n\right\}  }%
}\underbrace{\sum_{q=1}^{n}\left(  -1\right)  ^{p+q}c_{p,q}\det\left(  A_{\sim
p,\bullet}\right)  \det\left(  C_{\sim v,\sim q}\right)  }_{=\det\left(
A_{\sim p,\bullet}\right)  \sum_{q=1}^{n}\left(  -1\right)  ^{p+q}c_{p,q}%
\det\left(  C_{\sim v,\sim q}\right)  }\\
&  =\sum_{p\in\left\{  1,2,\ldots,n\right\}  }\det\left(  A_{\sim p,\bullet
}\right)  \sum_{q=1}^{n}\left(  -1\right)  ^{p+q}c_{p,q}\det\left(  C_{\sim
v,\sim q}\right) \\
&  =\det\left(  A_{\sim v,\bullet}\right)  \underbrace{\sum_{q=1}^{n}\left(
-1\right)  ^{v+q}c_{v,q}\det\left(  C_{\sim v,\sim q}\right)  }%
_{\substack{=\det C\\\text{(by (\ref{pf.prop.desnanot.AC.4}))}}}\\
&  \ \ \ \ \ \ \ \ \ \ +\sum_{\substack{p\in\left\{  1,2,\ldots,n\right\}
;\\p\neq v}}\det\left(  A_{\sim p,\bullet}\right)  \underbrace{\sum_{q=1}%
^{n}\left(  -1\right)  ^{p+q}c_{p,q}\det\left(  C_{\sim v,\sim q}\right)
}_{\substack{=0\\\text{(by (\ref{pf.prop.desnanot.AC.3}))}}}\\
&  \ \ \ \ \ \ \ \ \ \ \left(
\begin{array}
[c]{c}%
\text{here, we have split off the addend for }p=v\text{ from the sum,}\\
\text{since }v\in\left\{  1,2,\ldots,n\right\}
\end{array}
\right) \\
&  =\det\left(  A_{\sim v,\bullet}\right)  \det C+\underbrace{\sum
_{\substack{p\in\left\{  1,2,\ldots,n\right\}  ;\\p\neq v}}\det\left(  A_{\sim
p,\bullet}\right)  0}_{=0}=\det\left(  A_{\sim v,\bullet}\right)  \det C.
\end{align*}
This proves Proposition \ref{prop.desnanot.AC}.
\end{proof}

Now, let us show a purely technical lemma (gathering a few equalities for easy
access in a proof further down):

\begin{lemma}
\label{lem.desnanot.AB.tech}Let $n$ be a positive integer. Let $B\in
\mathbb{K}^{n\times\left(  n-1\right)  }$. Let $u$ and $v$ be two elements of
$\left\{  1,2,\ldots,n\right\}  $ such that $u<v$.

Consider the vector $\left(  I_{n}\right)  _{\bullet,u}\in\mathbb{K}%
^{n\times1}$. (This is the $u$-th column of the identity matrix $I_{n}$.\ \ \ \ \footnotemark)

Define an $n\times n$-matrix $C\in\mathbb{K}^{n\times n}$ by
\[
C=\left(  B\mid\left(  I_{n}\right)  _{\bullet,u}\right)  .
\]
Then, the following holds:

\textbf{(a)} We have%
\[
\det\left(  C_{\sim v,\sim q}\right)  =-\left(  -1\right)  ^{n+u}\det\left(
\operatorname*{rows}\nolimits_{1,2,\ldots,\widehat{u},\ldots,\widehat{v}%
,\ldots,n}\left(  B_{\bullet,\sim q}\right)  \right)
\]
for every $q\in\left\{  1,2,\ldots,n-1\right\}  $.

\textbf{(b)} We have%
\begin{equation}
\left(  -1\right)  ^{n+q}\det\left(  C_{\sim v,\sim q}\right)  =-\left(
-1\right)  ^{q+u}\det\left(  \operatorname*{rows}\nolimits_{1,2,\ldots
,\widehat{u},\ldots,\widehat{v},\ldots,n}\left(  B_{\bullet,\sim q}\right)
\right)  \label{pf.prop.desnanot.AB.1}%
\end{equation}
for every $q\in\left\{  1,2,\ldots,n-1\right\}  $.

\textbf{(c)} We have%
\begin{equation}
C_{\sim v,\sim n}=B_{\sim v,\bullet}. \label{pf.prop.desnanot.AB.2}%
\end{equation}

\textbf{(d)} We have
\begin{equation}
C_{\bullet,q}=B_{\bullet,q} \label{pf.prop.desnanot.AB.3}%
\end{equation}
for every $q\in\left\{  1,2,\ldots,n-1\right\}  $.

\textbf{(e)} Any $A\in\mathbb{K}^{n\times\left(  n-1\right)  }$ satisfies%
\begin{equation}
\det\left(  A\mid C_{\bullet,n}\right)  =\left(  -1\right)  ^{n+u}\det\left(
A_{\sim u,\bullet}\right)  . \label{pf.prop.desnanot.AB.4}%
\end{equation}

\textbf{(f)} We have
\begin{equation}
\det C=\left(  -1\right)  ^{n+u}\det\left(  B_{\sim u,\bullet}\right)  .
\label{pf.prop.desnanot.AB.5}%
\end{equation}

\end{lemma}

\footnotetext{Explicitly,
\[
\left(  I_{n}\right)  _{\bullet,u}=\left(  \underbrace{0,0,\ldots
,0}_{u-1\text{ zeroes}},1,\underbrace{0,0,\ldots,0}_{n-u\text{ zeroes}%
}\right)  ^{T}.
\]
}

\begin{proof}
[Proof of Lemma \ref{lem.desnanot.AB.tech}.]We have $n-1\in\mathbb{N}$ (since
$n$ is a positive integer). Also, $u<v\leq n$ (since $v\in\left\{
1,2,\ldots,n\right\}  $) and thus $u\leq n-1$ (since $u$ and $n$ are
integers). Combining this with $u\geq1$ (since $u\in\left\{  1,2,\ldots
,n\right\}  $), we obtain $u\in\left\{  1,2,\ldots,n-1\right\}  $.

\textbf{(a)} Let $q\in\left\{  1,2,\ldots,n-1\right\}  $. From $C=\left(
B\mid\left(  I_{n}\right)  _{\bullet,u}\right)  $, we obtain%
\[
C_{\bullet,\sim q}=\left(  B\mid\left(  I_{n}\right)  _{\bullet,u}\right)
_{\bullet,\sim q}=\left(  B_{\bullet,\sim q}\mid\left(  I_{n}\right)
_{\bullet,u}\right)
\]
(by Proposition \ref{prop.addcol.props1} \textbf{(c)}, applied to $n-1$, $B$
and $\left(  I_{n}\right)  _{\bullet,u}$ instead of $m$, $A$ and $v$). But
Proposition \ref{prop.unrows.basics} \textbf{(c)} (applied to $n$, $n$, $C$,
$v$ and $q$ instead of $n$, $m$, $A$, $u$ and $v$) yields $\left(
C_{\bullet,\sim q}\right)  _{\sim v,\bullet}=\left(  C_{\sim v,\bullet
}\right)  _{\bullet,\sim q}=C_{\sim v,\sim q}$. Hence,%
\begin{align}
C_{\sim v,\sim q}  &  =\left(  \underbrace{C_{\bullet,\sim q}}_{=\left(
B_{\bullet,\sim q}\mid\left(  I_{n}\right)  _{\bullet,u}\right)  }\right)
_{\sim v,\bullet}=\left(  B_{\bullet,\sim q}\mid\left(  I_{n}\right)
_{\bullet,u}\right)  _{\sim v,\bullet}\nonumber\\
&  =\left(  \left(  B_{\bullet,\sim q}\right)  _{\sim v,\bullet}\mid\left(
\left(  I_{n}\right)  _{\bullet,u}\right)  _{\sim v,\bullet}\right)
\label{pf.prop.desnanot.AB.1.pf.1}%
\end{align}
(by Proposition \ref{prop.addcol.props1} \textbf{(e)}, applied to $n$, $n-2$,
$B_{\bullet,\sim q}$, $\left(  I_{n}\right)  _{\bullet,u}$ and $v$ instead of
$n$, $m$, $A$, $v$ and $p$).

On the other hand, Proposition \ref{prop.unrows.basics} \textbf{(c)} (applied
to $n$, $n-1$, $B$, $v$ and $q$ instead of $n$, $m$, $A$, $u$ and $v$) yields
$\left(  B_{\bullet,\sim q}\right)  _{\sim v,\bullet}=\left(  B_{\sim
v,\bullet}\right)  _{\bullet,\sim q}=B_{\sim v,\sim q}$. Now,
(\ref{pf.prop.desnanot.AB.1.pf.1}) becomes%
\[
C_{\sim v,\sim q}=\left(  \underbrace{\left(  B_{\bullet,\sim q}\right)
_{\sim v,\bullet}}_{=\left(  B_{\sim v,\bullet}\right)  _{\bullet,\sim q}}%
\mid\underbrace{\left(  \left(  I_{n}\right)  _{\bullet,u}\right)  _{\sim
v,\bullet}}_{\substack{=\left(  I_{n-1}\right)  _{\bullet,u}\\\text{(by
Proposition \ref{prop.unrows.basics-I})}}}\right)  =\left(  \left(  B_{\sim
v,\bullet}\right)  _{\bullet,\sim q}\mid\left(  I_{n-1}\right)  _{\bullet
,u}\right)  .
\]
Hence,%
\begin{align*}
\det\underbrace{\left(  C_{\sim v,\sim q}\right)  }_{=\left(  \left(  B_{\sim
v,\bullet}\right)  _{\bullet,\sim q}\mid\left(  I_{n-1}\right)  _{\bullet
,u}\right)  }  &  =\det\left(  \left(  B_{\sim v,\bullet}\right)
_{\bullet,\sim q}\mid\left(  I_{n-1}\right)  _{\bullet,u}\right) \\
&  =\left(  -1\right)  ^{\left(  n-1\right)  +u}\det\left(  \left(  B_{\sim
v,\bullet}\right)  _{\sim u,\sim q}\right)
\end{align*}
(by Proposition \ref{prop.addcol.props3} \textbf{(d)}, applied to $n-1$,
$B_{\sim v,\bullet}$ and $u$ instead of $n$, $A$ and $p$). Thus,%
\begin{align*}
\det\left(  C_{\sim v,\sim q}\right)   &  =\underbrace{\left(  -1\right)
^{\left(  n-1\right)  +u}}_{\substack{=\left(  -1\right)  ^{n+u+1}%
\\\text{(since }\left(  n-1\right)  +u=\left(  n+u+1\right)  -2\\\equiv
n+u+1\operatorname{mod}2\text{)}}}\det\left(  \underbrace{\left(  B_{\sim
v,\bullet}\right)  _{\sim u,\sim q}}_{\substack{=\operatorname*{rows}%
\nolimits_{1,2,\ldots,\widehat{u},\ldots,\widehat{v},\ldots,n}\left(
B_{\bullet,\sim q}\right)  \\\text{(by Proposition \ref{prop.unrows.basics}
\textbf{(l)},}\\\text{applied to }B\text{ and }n-1\text{ instead of }A\text{
and }m\text{)}}}\right) \\
&  =\underbrace{\left(  -1\right)  ^{n+u+1}}_{=-\left(  -1\right)  ^{n+u}}%
\det\left(  \operatorname*{rows}\nolimits_{1,2,\ldots,\widehat{u}%
,\ldots,\widehat{v},\ldots,n}\left(  B_{\bullet,\sim q}\right)  \right) \\
&  =-\left(  -1\right)  ^{n+u}\det\left(  \operatorname*{rows}%
\nolimits_{1,2,\ldots,\widehat{u},\ldots,\widehat{v},\ldots,n}\left(
B_{\bullet,\sim q}\right)  \right)  .
\end{align*}
This proves Lemma \ref{lem.desnanot.AB.tech} \textbf{(a)}.

\textbf{(b)} Let $q\in\left\{  1,2,\ldots,n-1\right\}  $. Then,%
\begin{align*}
&  \left(  -1\right)  ^{n+q}\underbrace{\det\left(  C_{\sim v,\sim q}\right)
}_{\substack{=-\left(  -1\right)  ^{n+u}\det\left(  \operatorname*{rows}%
\nolimits_{1,2,\ldots,\widehat{u},\ldots,\widehat{v},\ldots,n}\left(
B_{\bullet,\sim q}\right)  \right)  \\\text{(by Lemma
\ref{lem.desnanot.AB.tech} \textbf{(a)})}}}\\
&  =\left(  -1\right)  ^{n+q}\left(  -\left(  -1\right)  ^{n+u}\det\left(
\operatorname*{rows}\nolimits_{1,2,\ldots,\widehat{u},\ldots,\widehat{v}%
,\ldots,n}\left(  B_{\bullet,\sim q}\right)  \right)  \right) \\
&  =-\underbrace{\left(  -1\right)  ^{n+q}\left(  -1\right)  ^{n+u}%
}_{\substack{=\left(  -1\right)  ^{\left(  n+q\right)  +\left(  n+u\right)
}=\left(  -1\right)  ^{q+u}\\\text{(since }\left(  n+q\right)  +\left(
n+u\right)  =2n+q+u\equiv q+u\operatorname{mod}2\text{)}}}\det\left(
\operatorname*{rows}\nolimits_{1,2,\ldots,\widehat{u},\ldots,\widehat{v}%
,\ldots,n}\left(  B_{\bullet,\sim q}\right)  \right) \\
&  =-\left(  -1\right)  ^{q+u}\det\left(  \operatorname*{rows}%
\nolimits_{1,2,\ldots,\widehat{u},\ldots,\widehat{v},\ldots,n}\left(
B_{\bullet,\sim q}\right)  \right)  .
\end{align*}
This proves Lemma \ref{lem.desnanot.AB.tech} \textbf{(b)}.

\textbf{(c)} We have $C=\left(  B\mid\left(  I_{n}\right)  _{\bullet
,u}\right)  $. Thus,%
\begin{align*}
C_{\sim v,\sim n}  &  =\left(  B\mid\left(  I_{n}\right)  _{\bullet,u}\right)
_{\sim v,\sim n}=\left(  B\mid\left(  I_{n}\right)  _{\bullet,u}\right)
_{\sim v,\sim\left(  \left(  n-1\right)  +1\right)  }%
\ \ \ \ \ \ \ \ \ \ \left(  \text{since }n=\left(  n-1\right)  +1\right) \\
&  =B_{\sim v,\bullet}\ \ \ \ \ \ \ \ \ \ \left(
\begin{array}
[c]{c}%
\text{by Proposition \ref{prop.addcol.props1} \textbf{(f)}, applied to}\\
n-1\text{, }B\text{ and }\left(  I_{n}\right)  _{\bullet,u}\text{ instead of
}m\text{, }A\text{ and }v
\end{array}
\right)  .
\end{align*}
This proves Lemma \ref{lem.desnanot.AB.tech} \textbf{(c)}.

\textbf{(d)} Let $q\in\left\{  1,2,\ldots,n-1\right\}  $. From $C=\left(
B\mid\left(  I_{n}\right)  _{\bullet,u}\right)  $, we obtain%
\[
C_{\bullet,q}=\left(  B\mid\left(  I_{n}\right)  _{\bullet,u}\right)
_{\bullet,q}=B_{\bullet,q}%
\]
(by Proposition \ref{prop.addcol.props1} \textbf{(a)}, applied to $n-1$, $B$
and $\left(  I_{n}\right)  _{\bullet,u}$ instead of $m$, $A$ and $v$). This
proves Lemma \ref{lem.desnanot.AB.tech} \textbf{(d)}.

\textbf{(e)} Let $A\in\mathbb{K}^{\left(  n-1\right)  \times n}$. Proposition
\ref{prop.addcol.props1} \textbf{(b)} (applied to $n-1$, $B$ and $\left(
I_{n}\right)  _{\bullet,u}$ instead of $m$, $A$ and $v$) yields $\left(
B\mid\left(  I_{n}\right)  _{\bullet,u}\right)  _{\bullet,\left(  n-1\right)
+1}=\left(  I_{n}\right)  _{\bullet,u}$. This rewrites as $\left(
B\mid\left(  I_{n}\right)  _{\bullet,u}\right)  _{\bullet,n}=\left(
I_{n}\right)  _{\bullet,u}$ (since $\left(  n-1\right)  +1=n$). Now,
$C=\left(  B\mid\left(  I_{n}\right)  _{\bullet,u}\right)  $, so that%
\[
C_{\bullet,n}=\left(  B\mid\left(  I_{n}\right)  _{\bullet,u}\right)
_{\bullet,n}=\left(  I_{n}\right)  _{\bullet,u}.
\]
Hence,%
\[
\det\left(  A\mid\underbrace{C_{\bullet,n}}_{=\left(  I_{n}\right)
_{\bullet,u}}\right)  =\det\left(  A\mid\left(  I_{n}\right)  _{\bullet
,u}\right)  =\left(  -1\right)  ^{n+u}\det\left(  A_{\sim u,\bullet}\right)
\]
(by Proposition \ref{prop.addcol.props2} \textbf{(b)}, applied to $p=u$). This
proves Lemma \ref{lem.desnanot.AB.tech} \textbf{(e)}.

\textbf{(f)} We have%
\[
\det\underbrace{C}_{=\left(  B\mid\left(  I_{n}\right)  _{\bullet,u}\right)
}=\det\left(  B\mid\left(  I_{n}\right)  _{\bullet,u}\right)  =\left(
-1\right)  ^{n+u}\det\left(  B_{\sim u,\bullet}\right)
\]
(by Proposition \ref{prop.addcol.props2} \textbf{(b)}, applied to $B$ and $u$
instead of $A$ and $p$). This proves Lemma \ref{lem.desnanot.AB.tech}
\textbf{(f)}.
\end{proof}

Next, we claim the following:

\begin{proposition}
\label{prop.desnanot.AB}Let $n$ be a positive integer. Let $A\in
\mathbb{K}^{n\times\left(  n-1\right)  }$ and $B\in\mathbb{K}^{n\times\left(
n-1\right)  }$. Let $u$ and $v$ be two elements of $\left\{  1,2,\ldots
,n\right\}  $ such that $u<v$. Then,%
\begin{align*}
&  \sum_{r=1}^{n-1}\left(  -1\right)  ^{r}\det\left(  A\mid B_{\bullet
,r}\right)  \det\left(  \operatorname*{rows}\nolimits_{1,2,\ldots
,\widehat{u},\ldots,\widehat{v},\ldots,n}\left(  B_{\bullet,\sim r}\right)
\right) \\
&  =\left(  -1\right)  ^{n}\left(  \det\left(  A_{\sim u,\bullet}\right)
\det\left(  B_{\sim v,\bullet}\right)  -\det\left(  A_{\sim v,\bullet}\right)
\det\left(  B_{\sim u,\bullet}\right)  \right)  .
\end{align*}

\end{proposition}

\begin{example}
\label{exam.prop.desnanot.AB}If we set $n=3$, $A=\left(
\begin{array}
[c]{cc}%
a & a^{\prime}\\
b & b^{\prime}\\
c & c^{\prime}%
\end{array}
\right)  $, $B=\left(
\begin{array}
[c]{cc}%
x & x^{\prime}\\
y & y^{\prime}\\
z & z^{\prime}%
\end{array}
\right)  $, $u=1$ and $v=2$, then Proposition \ref{prop.desnanot.AB} claims
that%
\begin{align*}
&  -\det\left(
\begin{array}
[c]{ccc}%
a & a^{\prime} & x\\
b & b^{\prime} & y\\
c & c^{\prime} & z
\end{array}
\right)  \det\left(
\begin{array}
[c]{c}%
z^{\prime}%
\end{array}
\right)  +\det\left(
\begin{array}
[c]{ccc}%
a & a^{\prime} & x^{\prime}\\
b & b^{\prime} & y^{\prime}\\
c & c^{\prime} & z^{\prime}%
\end{array}
\right)  \det\left(
\begin{array}
[c]{c}%
z
\end{array}
\right) \\
&  =\left(  -1\right)  ^{3}\left(  \det\left(
\begin{array}
[c]{cc}%
b & b^{\prime}\\
c & c^{\prime}%
\end{array}
\right)  \det\left(
\begin{array}
[c]{cc}%
x & x^{\prime}\\
z & z^{\prime}%
\end{array}
\right)  -\det\left(
\begin{array}
[c]{cc}%
a & a^{\prime}\\
c & c^{\prime}%
\end{array}
\right)  \det\left(
\begin{array}
[c]{cc}%
y & y^{\prime}\\
z & z^{\prime}%
\end{array}
\right)  \right)  .
\end{align*}

\end{example}

\begin{proof}
[Proof of Proposition \ref{prop.desnanot.AB}.]We have $n-1\in\mathbb{N}$
(since $n$ is a positive integer).

Define an $n\times n$-matrix $C\in\mathbb{K}^{n\times n}$ as in Lemma
\ref{lem.desnanot.AB.tech}.

Proposition \ref{prop.desnanot.AC} yields%
\begin{align*}
&  \det\left(  A_{\sim v,\bullet}\right)  \det C\\
&  =\sum_{q=1}^{n}\underbrace{\left(  -1\right)  ^{n+q}\det\left(  A\mid
C_{\bullet,q}\right)  }_{=\det\left(  A\mid C_{\bullet,q}\right)  \left(
-1\right)  ^{n+q}}\det\left(  C_{\sim v,\sim q}\right)  =\sum_{q=1}^{n}%
\det\left(  A\mid C_{\bullet,q}\right)  \left(  -1\right)  ^{n+q}\det\left(
C_{\sim v,\sim q}\right) \\
&  =\sum_{q=1}^{n-1}\det\left(  A\mid\underbrace{C_{\bullet,q}}%
_{\substack{=B_{\bullet,q}\\\text{(by (\ref{pf.prop.desnanot.AB.3}))}%
}}\right)  \underbrace{\left(  -1\right)  ^{n+q}\det\left(  C_{\sim v,\sim
q}\right)  }_{\substack{=-\left(  -1\right)  ^{q+u}\det\left(
\operatorname*{rows}\nolimits_{1,2,\ldots,\widehat{u},\ldots,\widehat{v}%
,\ldots,n}\left(  B_{\bullet,\sim q}\right)  \right)  \\\text{(by
(\ref{pf.prop.desnanot.AB.1}))}}}\\
&  \ \ \ \ \ \ \ \ \ \ +\underbrace{\det\left(  A\mid C_{\bullet,n}\right)
}_{\substack{=\left(  -1\right)  ^{n+u}\det\left(  A_{\sim u,\bullet}\right)
\\\text{(by (\ref{pf.prop.desnanot.AB.4}))}}}\underbrace{\left(  -1\right)
^{n+n}}_{\substack{=1\\\text{(since }n+n=2n\text{ is even)}}}\det\left(
\underbrace{C_{\sim v,\sim n}}_{\substack{=B_{\sim v,\bullet}\\\text{(by
(\ref{pf.prop.desnanot.AB.2}))}}}\right) \\
&  \ \ \ \ \ \ \ \ \ \ \left(  \text{here, we have split off the addend for
}q=n\text{ from the sum}\right) \\
&  =\underbrace{\sum_{q=1}^{n-1}\det\left(  A\mid B_{\bullet,q}\right)
\left(  -\left(  -1\right)  ^{q+u}\det\left(  \operatorname*{rows}%
\nolimits_{1,2,\ldots,\widehat{u},\ldots,\widehat{v},\ldots,n}\left(
B_{\bullet,\sim q}\right)  \right)  \right)  }_{=-\sum_{q=1}^{n-1}\left(
-1\right)  ^{q+u}\det\left(  A\mid B_{\bullet,q}\right)  \det\left(
\operatorname*{rows}\nolimits_{1,2,\ldots,\widehat{u},\ldots,\widehat{v}%
,\ldots,n}\left(  B_{\bullet,\sim q}\right)  \right)  }\\
&  \ \ \ \ \ \ \ \ \ \ +\left(  -1\right)  ^{n+u}\det\left(  A_{\sim
u,\bullet}\right)  \det\left(  B_{\sim v,\bullet}\right) \\
&  =-\sum_{q=1}^{n-1}\left(  -1\right)  ^{q+u}\det\left(  A\mid B_{\bullet
,q}\right)  \det\left(  \operatorname*{rows}\nolimits_{1,2,\ldots
,\widehat{u},\ldots,\widehat{v},\ldots,n}\left(  B_{\bullet,\sim q}\right)
\right) \\
&  \ \ \ \ \ \ \ \ \ \ +\left(  -1\right)  ^{n+u}\det\left(  A_{\sim
u,\bullet}\right)  \det\left(  B_{\sim v,\bullet}\right)  .
\end{align*}
Adding $\sum_{q=1}^{n-1}\left(  -1\right)  ^{q+u}\det\left(  A\mid
B_{\bullet,q}\right)  \det\left(  \operatorname*{rows}\nolimits_{1,2,\ldots
,\widehat{u},\ldots,\widehat{v},\ldots,n}\left(  B_{\bullet,\sim q}\right)
\right)  $ to both sides of this equality, we obtain%
\begin{align*}
&  \sum_{q=1}^{n-1}\left(  -1\right)  ^{q+u}\det\left(  A\mid B_{\bullet
,q}\right)  \det\left(  \operatorname*{rows}\nolimits_{1,2,\ldots
,\widehat{u},\ldots,\widehat{v},\ldots,n}\left(  B_{\bullet,\sim q}\right)
\right)  +\det\left(  A_{\sim v,\bullet}\right)  \det C\\
&  =\left(  -1\right)  ^{n+u}\det\left(  A_{\sim u,\bullet}\right)
\det\left(  B_{\sim v,\bullet}\right)  .
\end{align*}
Subtracting $\det\left(  A_{\sim v,\bullet}\right)  \det C$ from both sides of
this equality, we find%
\begin{align*}
&  \sum_{q=1}^{n-1}\left(  -1\right)  ^{q+u}\det\left(  A\mid B_{\bullet
,q}\right)  \det\left(  \operatorname*{rows}\nolimits_{1,2,\ldots
,\widehat{u},\ldots,\widehat{v},\ldots,n}\left(  B_{\bullet,\sim q}\right)
\right) \\
&  =\left(  -1\right)  ^{n+u}\det\left(  A_{\sim u,\bullet}\right)
\det\left(  B_{\sim v,\bullet}\right)  -\det\left(  A_{\sim v,\bullet}\right)
\underbrace{\det C}_{\substack{=\left(  -1\right)  ^{n+u}\det\left(  B_{\sim
u,\bullet}\right)  \\\text{(by (\ref{pf.prop.desnanot.AB.5}))}}}\\
&  =\left(  -1\right)  ^{n+u}\det\left(  A_{\sim u,\bullet}\right)
\det\left(  B_{\sim v,\bullet}\right)  -\underbrace{\det\left(  A_{\sim
v,\bullet}\right)  \left(  -1\right)  ^{n+u}}_{=\left(  -1\right)  ^{n+u}%
\det\left(  A_{\sim v,\bullet}\right)  }\det\left(  B_{\sim u,\bullet}\right)
\\
&  =\left(  -1\right)  ^{n+u}\det\left(  A_{\sim u,\bullet}\right)
\det\left(  B_{\sim v,\bullet}\right)  -\left(  -1\right)  ^{n+u}\det\left(
A_{\sim v,\bullet}\right)  \det\left(  B_{\sim u,\bullet}\right) \\
&  =\left(  -1\right)  ^{n+u}\left(  \det\left(  A_{\sim u,\bullet}\right)
\det\left(  B_{\sim v,\bullet}\right)  -\det\left(  A_{\sim v,\bullet}\right)
\det\left(  B_{\sim u,\bullet}\right)  \right)  .
\end{align*}
Multiplying both sides of this equality by $\left(  -1\right)  ^{u}$, we
obtain%
\begin{align*}
&  \left(  -1\right)  ^{u}\sum_{q=1}^{n-1}\left(  -1\right)  ^{q+u}\det\left(
A\mid B_{\bullet,q}\right)  \det\left(  \operatorname*{rows}%
\nolimits_{1,2,\ldots,\widehat{u},\ldots,\widehat{v},\ldots,n}\left(
B_{\bullet,\sim q}\right)  \right) \\
&  =\underbrace{\left(  -1\right)  ^{u}\left(  -1\right)  ^{n+u}%
}_{\substack{=\left(  -1\right)  ^{u+\left(  n+u\right)  }=\left(  -1\right)
^{n}\\\text{(since }u+\left(  n+u\right)  =2u+n\equiv n\operatorname{mod}%
2\text{)}}}\left(  \det\left(  A_{\sim u,\bullet}\right)  \det\left(  B_{\sim
v,\bullet}\right)  -\det\left(  A_{\sim v,\bullet}\right)  \det\left(  B_{\sim
u,\bullet}\right)  \right) \\
&  =\left(  -1\right)  ^{n}\left(  \det\left(  A_{\sim u,\bullet}\right)
\det\left(  B_{\sim v,\bullet}\right)  -\det\left(  A_{\sim v,\bullet}\right)
\det\left(  B_{\sim u,\bullet}\right)  \right)  ,
\end{align*}
so that%
\begin{align*}
&  \left(  -1\right)  ^{n}\left(  \det\left(  A_{\sim u,\bullet}\right)
\det\left(  B_{\sim v,\bullet}\right)  -\det\left(  A_{\sim v,\bullet}\right)
\det\left(  B_{\sim u,\bullet}\right)  \right) \\
&  =\left(  -1\right)  ^{u}\sum_{q=1}^{n-1}\left(  -1\right)  ^{q+u}%
\det\left(  A\mid B_{\bullet,q}\right)  \det\left(  \operatorname*{rows}%
\nolimits_{1,2,\ldots,\widehat{u},\ldots,\widehat{v},\ldots,n}\left(
B_{\bullet,\sim q}\right)  \right) \\
&  =\sum_{q=1}^{n-1}\underbrace{\left(  -1\right)  ^{u}\left(  -1\right)
^{q+u}}_{\substack{=\left(  -1\right)  ^{u+\left(  q+u\right)  }=\left(
-1\right)  ^{q}\\\text{(since }u+\left(  q+u\right)  =2u+q\equiv
q\operatorname{mod}2\text{)}}}\det\left(  A\mid B_{\bullet,q}\right)
\det\left(  \operatorname*{rows}\nolimits_{1,2,\ldots,\widehat{u}%
,\ldots,\widehat{v},\ldots,n}\left(  B_{\bullet,\sim q}\right)  \right) \\
&  =\sum_{q=1}^{n-1}\left(  -1\right)  ^{q}\det\left(  A\mid B_{\bullet
,q}\right)  \det\left(  \operatorname*{rows}\nolimits_{1,2,\ldots
,\widehat{u},\ldots,\widehat{v},\ldots,n}\left(  B_{\bullet,\sim q}\right)
\right) \\
&  =\sum_{r=1}^{n-1}\left(  -1\right)  ^{r}\det\left(  A\mid B_{\bullet
,r}\right)  \det\left(  \operatorname*{rows}\nolimits_{1,2,\ldots
,\widehat{u},\ldots,\widehat{v},\ldots,n}\left(  B_{\bullet,\sim r}\right)
\right)
\end{align*}
(here, we have renamed the summation index $q$ as $r$). This proves
Proposition \ref{prop.desnanot.AB}.
\end{proof}

Now, we can finally prove Theorem \ref{thm.desnanot}:

\begin{proof}
[Proof of Theorem \ref{thm.desnanot}.]We have $v\in\left\{  1,2,\ldots
,n\right\}  $ and thus $v\leq n$. Hence, $u<v\leq n$, so that $u\leq n-1$
(since $u$ and $n$ are integers). Also, $u\in\left\{  1,2,\ldots,n\right\}  $,
so that $1\leq u$. Combining $1\leq u$ with $u\leq n-1$, we obtain
$u\in\left\{  1,2,\ldots,n-1\right\}  $.

Proposition \ref{prop.submatrix.easy} \textbf{(d)} (applied to $n$, $n-2$,
$\left(  1,2,\ldots,\widehat{p},\ldots,\widehat{q},\ldots,n\right)  $, $n-2$
and $\left(  1,2,\ldots,\widehat{u},\ldots,\widehat{v},\ldots,n\right)  $
instead of $m$, $u$, $\left(  i_{1},i_{2},\ldots,i_{u}\right)  $, $v$ and
$\left(  j_{1},j_{2},\ldots,j_{v}\right)  $) yields%
\begin{align}
\operatorname*{sub}\nolimits_{1,2,\ldots,\widehat{p},\ldots,\widehat{q}%
,\ldots,n}^{1,2,\ldots,\widehat{u},\ldots,\widehat{v},\ldots,n}A  &
=\operatorname*{rows}\nolimits_{1,2,\ldots,\widehat{p},\ldots,\widehat{q}%
,\ldots,n}\left(  \operatorname*{cols}\nolimits_{1,2,\ldots,\widehat{u}%
,\ldots,\widehat{v},\ldots,n}A\right) \label{pf.thm.desnanot.1}\\
&  =\operatorname*{cols}\nolimits_{1,2,\ldots,\widehat{u},\ldots
,\widehat{v},\ldots,n}\left(  \operatorname*{rows}\nolimits_{1,2,\ldots
,\widehat{p},\ldots,\widehat{q},\ldots,n}A\right)  .\nonumber
\end{align}

We have $u<v$, so that $u\leq v-1$ (since $u$ and $v$ are integers). Combining
this with $1\leq u$, we obtain $u\in\left\{  1,2,\ldots,v-1\right\}  $. Thus,
Proposition \ref{prop.unrows.basics} \textbf{(h)} (applied to $m=n$ and $w=u$)
yields $\left(  A_{\bullet,\sim v}\right)  _{\bullet,\sim u}%
=\operatorname*{cols}\nolimits_{1,2,\ldots,\widehat{u},\ldots,\widehat{v}%
,\ldots,n}A$. Thus,%
\begin{align}
&  \operatorname*{rows}\nolimits_{1,2,\ldots,\widehat{p},\ldots,\widehat{q}%
,\ldots,n}\left(  \underbrace{\left(  A_{\bullet,\sim v}\right)
_{\bullet,\sim u}}_{=\operatorname*{cols}\nolimits_{1,2,\ldots,\widehat{u}%
,\ldots,\widehat{v},\ldots,n}A}\right) \nonumber\\
&  =\operatorname*{rows}\nolimits_{1,2,\ldots,\widehat{p},\ldots
,\widehat{q},\ldots,n}\left(  \operatorname*{cols}\nolimits_{1,2,\ldots
,\widehat{u},\ldots,\widehat{v},\ldots,n}A\right) \nonumber\\
&  =\operatorname*{sub}\nolimits_{1,2,\ldots,\widehat{p},\ldots,\widehat{q}%
,\ldots,n}^{1,2,\ldots,\widehat{u},\ldots,\widehat{v},\ldots,n}%
A\ \ \ \ \ \ \ \ \ \ \left(  \text{by (\ref{pf.thm.desnanot.1})}\right)  .
\label{pf.thm.desnanot.2}%
\end{align}

\begin{vershort}
Proposition \ref{prop.unrows.basics} \textbf{(c)} (applied to $n$, $p$ and $u$
instead of $m$, $u$ and $v$) yields $\left(  A_{\bullet,\sim u}\right)  _{\sim
p,\bullet}=\left(  A_{\sim p,\bullet}\right)  _{\bullet,\sim u}=A_{\sim p,\sim
u}$. Similarly, $\left(  A_{\bullet,\sim v}\right)  _{\sim p,\bullet}=\left(
A_{\sim p,\bullet}\right)  _{\bullet,\sim v}=A_{\sim p,\sim v}$ and $\left(
A_{\bullet,\sim u}\right)  _{\sim q,\bullet}=\left(  A_{\sim q,\bullet
}\right)  _{\bullet,\sim u}=A_{\sim q,\sim u}$ and $\left(  A_{\bullet,\sim
v}\right)  _{\sim q,\bullet}=\left(  A_{\sim q,\bullet}\right)  _{\bullet,\sim
v}=A_{\sim q,\sim v}$.
\end{vershort}

\begin{verlong}
Proposition \ref{prop.unrows.basics} \textbf{(c)} (applied to $n$, $p$ and $u$
instead of $m$, $u$ and $v$) yields $\left(  A_{\bullet,\sim u}\right)  _{\sim
p,\bullet}=\left(  A_{\sim p,\bullet}\right)  _{\bullet,\sim u}=A_{\sim p,\sim
u}$.

Proposition \ref{prop.unrows.basics} \textbf{(c)} (applied to $n$, $p$ and $v$
instead of $m$, $u$ and $v$) yields $\left(  A_{\bullet,\sim v}\right)  _{\sim
p,\bullet}=\left(  A_{\sim p,\bullet}\right)  _{\bullet,\sim v}=A_{\sim p,\sim
v}$.

Proposition \ref{prop.unrows.basics} \textbf{(c)} (applied to $n$, $q$ and $u$
instead of $m$, $u$ and $v$) yields $\left(  A_{\bullet,\sim u}\right)  _{\sim
q,\bullet}=\left(  A_{\sim q,\bullet}\right)  _{\bullet,\sim u}=A_{\sim q,\sim
u}$.

Proposition \ref{prop.unrows.basics} \textbf{(c)} (applied to $n$, $q$ and $v$
instead of $m$, $u$ and $v$) yields $\left(  A_{\bullet,\sim v}\right)  _{\sim
q,\bullet}=\left(  A_{\sim q,\bullet}\right)  _{\bullet,\sim v}=A_{\sim q,\sim
v}$.
\end{verlong}

The integer $n$ is positive (since $n\geq2$). Thus, Proposition
\ref{prop.desnanot.AB} (applied to $A_{\bullet,\sim u}$, $A_{\bullet,\sim v}$,
$p$ and $q$ instead of $A$, $B$, $u$ and $v$) yields%
\begin{align*}
&  \sum_{r=1}^{n-1}\left(  -1\right)  ^{r}\det\left(  A_{\bullet,\sim u}%
\mid\left(  A_{\bullet,\sim v}\right)  _{\bullet,r}\right)  \det\left(
\operatorname*{rows}\nolimits_{1,2,\ldots,\widehat{p},\ldots,\widehat{q}%
,\ldots,n}\left(  \left(  A_{\bullet,\sim v}\right)  _{\bullet,\sim r}\right)
\right) \\
&  =\left(  -1\right)  ^{n}\left(  \det\left(  \underbrace{\left(
A_{\bullet,\sim u}\right)  _{\sim p,\bullet}}_{=A_{\sim p,\sim u}}\right)
\det\left(  \underbrace{\left(  A_{\bullet,\sim v}\right)  _{\sim q,\bullet}%
}_{=A_{\sim q,\sim v}}\right)  \right. \\
&  \ \ \ \ \ \ \ \ \ \ \ \ \ \ \ \ \ \ \ \ \left.  -\det\left(
\underbrace{\left(  A_{\bullet,\sim u}\right)  _{\sim q,\bullet}}_{=A_{\sim
q,\sim u}}\right)  \det\left(  \underbrace{\left(  A_{\bullet,\sim v}\right)
_{\sim p,\bullet}}_{=A_{\sim p,\sim v}}\right)  \right) \\
&  =\left(  -1\right)  ^{n}\left(  \det\left(  A_{\sim p,\sim u}\right)
\det\left(  A_{\sim q,\sim v}\right)  -\det\left(  A_{\sim q,\sim u}\right)
\det\left(  A_{\sim p,\sim v}\right)  \right)  .
\end{align*}
Hence,%
\begin{align*}
&  \left(  -1\right)  ^{n}\left(  \det\left(  A_{\sim p,\sim u}\right)
\det\left(  A_{\sim q,\sim v}\right)  -\det\left(  A_{\sim q,\sim u}\right)
\det\left(  A_{\sim p,\sim v}\right)  \right) \\
&  =\underbrace{\sum_{r=1}^{n-1}}_{=\sum_{r\in\left\{  1,2,\ldots,n-1\right\}
}}\left(  -1\right)  ^{r}\det\left(  A_{\bullet,\sim u}\mid\left(
A_{\bullet,\sim v}\right)  _{\bullet,r}\right)  \det\left(
\operatorname*{rows}\nolimits_{1,2,\ldots,\widehat{p},\ldots,\widehat{q}%
,\ldots,n}\left(  \left(  A_{\bullet,\sim v}\right)  _{\bullet,\sim r}\right)
\right) \\
&  =\sum_{r\in\left\{  1,2,\ldots,n-1\right\}  }\left(  -1\right)  ^{r}%
\det\left(  A_{\bullet,\sim u}\mid\left(  A_{\bullet,\sim v}\right)
_{\bullet,r}\right)  \det\left(  \operatorname*{rows}\nolimits_{1,2,\ldots
,\widehat{p},\ldots,\widehat{q},\ldots,n}\left(  \left(  A_{\bullet,\sim
v}\right)  _{\bullet,\sim r}\right)  \right) \\
&  =\sum_{\substack{r\in\left\{  1,2,\ldots,n-1\right\}  ;\\r\neq u}}\left(
-1\right)  ^{r}\underbrace{\det\left(  A_{\bullet,\sim u}\mid\left(
A_{\bullet,\sim v}\right)  _{\bullet,r}\right)  }_{\substack{=0\\\text{(by
Proposition \ref{prop.addcol.props3} \textbf{(e)})}}}\det\left(
\operatorname*{rows}\nolimits_{1,2,\ldots,\widehat{p},\ldots,\widehat{q}%
,\ldots,n}\left(  \left(  A_{\bullet,\sim v}\right)  _{\bullet,\sim r}\right)
\right) \\
&  \ \ \ \ \ \ \ \ \ \ +\underbrace{\left(  -1\right)  ^{u}\det\left(
A_{\bullet,\sim u}\mid\left(  A_{\bullet,\sim v}\right)  _{\bullet,u}\right)
}_{\substack{=\left(  -1\right)  ^{n}\det A\\\text{(by Proposition
\ref{prop.addcol.props3} \textbf{(f)})}}}\det\left(
\underbrace{\operatorname*{rows}\nolimits_{1,2,\ldots,\widehat{p}%
,\ldots,\widehat{q},\ldots,n}\left(  \left(  A_{\bullet,\sim v}\right)
_{\bullet,\sim u}\right)  }_{\substack{=\operatorname*{sub}%
\nolimits_{1,2,\ldots,\widehat{p},\ldots,\widehat{q},\ldots,n}^{1,2,\ldots
,\widehat{u},\ldots,\widehat{v},\ldots,n}A\\\text{(by (\ref{pf.thm.desnanot.2}%
))}}}\right) \\
&  \ \ \ \ \ \ \ \ \ \ \left(
\begin{array}
[c]{c}%
\text{here, we have split off the addend for }r=u\text{ from the sum,}\\
\text{since }u\in\left\{  1,2,\ldots,n-1\right\}
\end{array}
\right) \\
&  =\underbrace{\sum_{\substack{r\in\left\{  1,2,\ldots,n-1\right\}  ;\\r\neq
u}}\left(  -1\right)  ^{r}0\det\left(  \operatorname*{rows}%
\nolimits_{1,2,\ldots,\widehat{p},\ldots,\widehat{q},\ldots,n}\left(  \left(
A_{\bullet,\sim v}\right)  _{\bullet,\sim r}\right)  \right)  }_{=0}\\
&  \ \ \ \ \ \ \ \ \ \ +\left(  -1\right)  ^{n}\det A\cdot\det\left(
\operatorname*{sub}\nolimits_{1,2,\ldots,\widehat{p},\ldots,\widehat{q}%
,\ldots,n}^{1,2,\ldots,\widehat{u},\ldots,\widehat{v},\ldots,n}A\right) \\
&  =\left(  -1\right)  ^{n}\det A\cdot\det\left(  \operatorname*{sub}%
\nolimits_{1,2,\ldots,\widehat{p},\ldots,\widehat{q},\ldots,n}^{1,2,\ldots
,\widehat{u},\ldots,\widehat{v},\ldots,n}A\right)  .
\end{align*}
Multiplying both sides of this equality by $\left(  -1\right)  ^{n}$, we
obtain%
\begin{align*}
&  \left(  -1\right)  ^{n}\left(  -1\right)  ^{n}\left(  \det\left(  A_{\sim
p,\sim u}\right)  \det\left(  A_{\sim q,\sim v}\right)  -\det\left(  A_{\sim
q,\sim u}\right)  \det\left(  A_{\sim p,\sim v}\right)  \right) \\
&  =\underbrace{\left(  -1\right)  ^{n}\left(  -1\right)  ^{n}}%
_{\substack{=\left(  -1\right)  ^{n+n}=1\\\text{(since }n+n=2n\text{ is
even)}}}\det A\cdot\det\left(  \operatorname*{sub}\nolimits_{1,2,\ldots
,\widehat{p},\ldots,\widehat{q},\ldots,n}^{1,2,\ldots,\widehat{u}%
,\ldots,\widehat{v},\ldots,n}A\right) \\
&  =\det A\cdot\det\left(  \operatorname*{sub}\nolimits_{1,2,\ldots
,\widehat{p},\ldots,\widehat{q},\ldots,n}^{1,2,\ldots,\widehat{u}%
,\ldots,\widehat{v},\ldots,n}A\right)  .
\end{align*}
Thus,%
\begin{align*}
&  \det A\cdot\det\left(  \operatorname*{sub}\nolimits_{1,2,\ldots
,\widehat{p},\ldots,\widehat{q},\ldots,n}^{1,2,\ldots,\widehat{u}%
,\ldots,\widehat{v},\ldots,n}A\right) \\
&  =\underbrace{\left(  -1\right)  ^{n}\left(  -1\right)  ^{n}}%
_{\substack{=\left(  -1\right)  ^{n+n}=1\\\text{(since }n+n=2n\text{ is
even)}}}\left(  \det\left(  A_{\sim p,\sim u}\right)  \det\left(  A_{\sim
q,\sim v}\right)  -\det\left(  A_{\sim q,\sim u}\right)  \det\left(  A_{\sim
p,\sim v}\right)  \right) \\
&  =\det\left(  A_{\sim p,\sim u}\right)  \det\left(  A_{\sim q,\sim
v}\right)  -\det\left(  A_{\sim q,\sim u}\right)  \det\left(  A_{\sim p,\sim
v}\right)  .
\end{align*}
This proves Theorem \ref{thm.desnanot}.
\end{proof}

Now that Theorem \ref{thm.desnanot} is proven, we conclude that Proposition
\ref{prop.desnanot.12} and Proposition \ref{prop.desnanot.1n} hold as well
(because in Exercise \ref{exe.unrows.basics} \textbf{(b)}, these two
propositions have been derived from Theorem \ref{thm.desnanot}).

\begin{exercise}
\label{exe.desnanot.jaw}Let $n$ be a positive integer. Let $B\in
\mathbb{K}^{n\times\left(  n-1\right)  }$. Fix $q\in\left\{  1,2,\ldots
,n-1\right\}  $. For every $x\in\left\{  1,2,\ldots,n\right\}  $, set
\[
\alpha_{x}=\det\left(  B_{\sim x,\bullet}\right)  \text{.}%
\]
For every two elements $x$ and $y$ of $\left\{  1,2,\ldots,n\right\}  $
satisfying $x<y$, set%
\[
\beta_{x,y}=\det\left(  \operatorname*{rows}\nolimits_{1,2,\ldots
,\widehat{x},\ldots,\widehat{y},\ldots,n}\left(  B_{\bullet,\sim q}\right)
\right)  .
\]
(Note that this depends on $q$, but we do not mention $q$ in the notation
because $q$ is fixed.)

Let $u$, $v$ and $w$ be three elements of $\left\{  1,2,\ldots,n\right\}  $
such that $u<v<w$. Thus, $\beta_{u,v}$, $\beta_{v,w}$ and $\beta_{u,w}$ are
well-defined elements of $\mathbb{K}$. Prove that%
\[
\alpha_{u}\beta_{v,w}+\alpha_{w}\beta_{u,v}=\alpha_{v}\beta_{u,w}.
\]

\end{exercise}

\begin{example}
\label{exam.exe.desnanot.jaw}If we set $n=4$, $B=\left(
\begin{array}
[c]{ccc}%
a & a^{\prime} & a^{\prime\prime}\\
b & b^{\prime} & b^{\prime\prime}\\
c & c^{\prime} & c^{\prime\prime}\\
d & d^{\prime} & d^{\prime\prime}%
\end{array}
\right)  $, $q=3$, $u=1$, $v=2$ and $w=3$, then Exercise
\ref{exe.desnanot.jaw} says that%
\begin{align*}
&  \det\left(
\begin{array}
[c]{ccc}%
b & b^{\prime} & b^{\prime\prime}\\
c & c^{\prime} & c^{\prime\prime}\\
d & d^{\prime} & d^{\prime\prime}%
\end{array}
\right)  \cdot\det\left(
\begin{array}
[c]{cc}%
a & a^{\prime}\\
d & d^{\prime}%
\end{array}
\right)  +\det\left(
\begin{array}
[c]{ccc}%
a & a^{\prime} & a^{\prime\prime}\\
b & b^{\prime} & b^{\prime\prime}\\
d & d^{\prime} & d^{\prime\prime}%
\end{array}
\right)  \cdot\det\left(
\begin{array}
[c]{cc}%
c & c^{\prime}\\
d & d^{\prime}%
\end{array}
\right) \\
&  =\det\left(
\begin{array}
[c]{ccc}%
a & a^{\prime} & a^{\prime\prime}\\
c & c^{\prime} & c^{\prime\prime}\\
d & d^{\prime} & d^{\prime\prime}%
\end{array}
\right)  \cdot\det\left(
\begin{array}
[c]{cc}%
b & b^{\prime}\\
d & d^{\prime}%
\end{array}
\right)  .
\end{align*}

\end{example}

\begin{remark}
Exercise \ref{exe.desnanot.jaw} appears in \cite[proof of Theorem 9]%
{KenWil14}, where it (or, rather, a certain transformation of determinant
expressions that relies on it) is called the \textquotedblleft jaw
move\textquotedblright.
\end{remark}

\begin{exercise}
\label{exe.desnanot.skew}Let $n\in\mathbb{N}$. Let $A$ be an alternating
$n\times n$-matrix. (See Definition \ref{def.altern} \textbf{(b)} for what
this means.) Let $S$ be any $n\times n$-matrix. Prove that each entry of the
matrix $\left(  \operatorname*{adj}S\right)  ^{T}A\left(  \operatorname*{adj}%
S\right)  $ is a multiple of $\det S$.
\end{exercise}

\subsection{The Pl\"{u}cker relation}

The following section is devoted to the \textit{Pl\"{u}cker relations}, or,
rather, one of the many things that tend to carry this name in the
literature\footnote{Most of the relevant literature, unfortunately, is not
very elementary, as the Pl\"{u}cker relations are at their most useful in the
algebraic geometry of the Grassmannian and of flag varieties
(\textquotedblleft Schubert calculus\textquotedblright). See \cite{KleLak72},
\cite[\S 3.4]{Jacobs10} and \cite[\S 9.1]{Fulton-Young} for expositions (all
three, however, well above the level of the present notes).}. The proofs will
be fairly short, since we did much of the necessary work in Section
\ref{sect.desnanot} already.

We shall use the notations of Definition \ref{def.unrows} throughout this section.

We begin with the following identity:

\begin{proposition}
\label{prop.pluecker.pre.row}Let $n$ be a positive integer. Let $B\in
\mathbb{K}^{n\times\left(  n-1\right)  }$. Then:

\textbf{(a)} We have%
\[
\sum_{r=1}^{n}\left(  -1\right)  ^{r}\det\left(  B_{\sim r,\bullet}\right)
B_{r,\bullet}=0_{1\times\left(  n-1\right)  }.
\]
(Recall that the product $\det\left(  B_{\sim r,\bullet}\right)  B_{r,\bullet
}$ in this equality is the product of the scalar $\det\left(  B_{\sim
r,\bullet}\right)  \in\mathbb{K}$ with the row vector $B_{r,\bullet}%
\in\mathbb{K}^{1\times\left(  n-1\right)  }$; as all such products, it is
computed entrywise, i.e., by the formula $\lambda\left(  a_{1},a_{2}%
,\ldots,a_{n-1}\right)  =\left(  \lambda a_{1},\lambda a_{2},\ldots,\lambda
a_{n-1}\right)  $.)

\textbf{(b)} Write the matrix $B$ in the form $B=\left(  b_{i,j}\right)
_{1\leq i\leq n,\ 1\leq j\leq n-1}$. Then,
\[
\sum_{r=1}^{n}\left(  -1\right)  ^{r}\det\left(  B_{\sim r,\bullet}\right)
b_{r,q}=0
\]
for every $q\in\left\{  1,2,\ldots,n-1\right\}  $.
\end{proposition}

Note that part \textbf{(a)} of Proposition \ref{prop.pluecker.pre.row} claims
an equality between two row vectors (indeed, $B_{r,\bullet}$ is a row vector
with $n-1$ entries for each $r\in\left\{  1,2,\ldots,n\right\}  $), whereas
part \textbf{(b)} claims an equality between two elements of $\mathbb{K}$ (for
each $q\in\left\{  1,2,\ldots,n-1\right\}  $). That said, the two parts are
essentially restatements of one another, and we will derive part \textbf{(a)}
from part \textbf{(b)} soon enough. Let us first illustrate Proposition
\ref{prop.pluecker.pre.row} on an example:

\begin{example}
\label{exam.prop.pluecker.pre.row}For this example, set $n=4$ and $B=\left(
\begin{array}
[c]{ccc}%
a & a^{\prime} & a^{\prime\prime}\\
b & b^{\prime} & b^{\prime\prime}\\
c & c^{\prime} & c^{\prime\prime}\\
d & d^{\prime} & d^{\prime\prime}%
\end{array}
\right)  $. Then, Proposition \ref{prop.pluecker.pre.row} \textbf{(a)} says
that%
\begin{align*}
&  -\det\left(
\begin{array}
[c]{ccc}%
b & b^{\prime} & b^{\prime\prime}\\
c & c^{\prime} & c^{\prime\prime}\\
d & d^{\prime} & d^{\prime\prime}%
\end{array}
\right)  \cdot\left(
\begin{array}
[c]{ccc}%
a & a^{\prime} & a^{\prime\prime}%
\end{array}
\right)  +\det\left(
\begin{array}
[c]{ccc}%
a & a^{\prime} & a^{\prime\prime}\\
c & c^{\prime} & c^{\prime\prime}\\
d & d^{\prime} & d^{\prime\prime}%
\end{array}
\right)  \cdot\left(
\begin{array}
[c]{ccc}%
b & b^{\prime} & b^{\prime\prime}%
\end{array}
\right) \\
&  \ \ \ \ \ \ \ \ \ \ -\det\left(
\begin{array}
[c]{ccc}%
a & a^{\prime} & a^{\prime\prime}\\
b & b^{\prime} & b^{\prime\prime}\\
d & d^{\prime} & d^{\prime\prime}%
\end{array}
\right)  \cdot\left(
\begin{array}
[c]{ccc}%
c & c^{\prime} & c^{\prime\prime}%
\end{array}
\right)  +\det\left(
\begin{array}
[c]{ccc}%
a & a^{\prime} & a^{\prime\prime}\\
b & b^{\prime} & b^{\prime\prime}\\
c & c^{\prime} & c^{\prime\prime}%
\end{array}
\right)  \cdot\left(
\begin{array}
[c]{ccc}%
d & d^{\prime} & d^{\prime\prime}%
\end{array}
\right) \\
&  =0_{1\times3}.
\end{align*}
Proposition \ref{prop.pluecker.pre.row} \textbf{(b)} (applied to $q=3$) yields%
\begin{align*}
&  -\det\left(
\begin{array}
[c]{ccc}%
b & b^{\prime} & b^{\prime\prime}\\
c & c^{\prime} & c^{\prime\prime}\\
d & d^{\prime} & d^{\prime\prime}%
\end{array}
\right)  \cdot a^{\prime\prime}+\det\left(
\begin{array}
[c]{ccc}%
a & a^{\prime} & a^{\prime\prime}\\
c & c^{\prime} & c^{\prime\prime}\\
d & d^{\prime} & d^{\prime\prime}%
\end{array}
\right)  \cdot b^{\prime\prime}\\
&  \ \ \ \ \ \ \ \ \ \ -\det\left(
\begin{array}
[c]{ccc}%
a & a^{\prime} & a^{\prime\prime}\\
b & b^{\prime} & b^{\prime\prime}\\
d & d^{\prime} & d^{\prime\prime}%
\end{array}
\right)  \cdot c^{\prime\prime}+\det\left(
\begin{array}
[c]{ccc}%
a & a^{\prime} & a^{\prime\prime}\\
b & b^{\prime} & b^{\prime\prime}\\
c & c^{\prime} & c^{\prime\prime}%
\end{array}
\right)  \cdot d^{\prime\prime}\\
&  =0.
\end{align*}

\end{example}

\begin{proof}
[Proof of Proposition \ref{prop.pluecker.pre.row}.]\textbf{(b)} Let
$q\in\left\{  1,2,\ldots,n-1\right\}  $.

\begin{vershort}
We shall use the notation introduced in Definition \ref{def.addcol}. We know
that $B$ is an $n\times\left(  n-1\right)  $-matrix (since $B\in
\mathbb{K}^{n\times\left(  n-1\right)  }$). Thus, $\left(  B\mid B_{\bullet
,q}\right)  $ is an $n\times n$-matrix. This $n\times n$-matrix $\left(  B\mid
B_{\bullet,q}\right)  $ is defined as the $n\times\left(  \left(  n-1\right)
+1\right)  $-matrix whose columns are $B_{\bullet,1},B_{\bullet,2}%
,\ldots,B_{\bullet,n-1},B_{\bullet,q}$; thus, it has two equal columns
(indeed, the column vector $B_{\bullet,q}$ appears twice among the columns
$B_{\bullet,1},B_{\bullet,2},\ldots,B_{\bullet,n-1},B_{\bullet,q}$). Thus,
Exercise \ref{exe.ps4.6} \textbf{(f)} (applied to $A=\left(  B\mid
B_{\bullet,q}\right)  $) shows that $\det\left(  B\mid B_{\bullet,q}\right)
=0$.
\end{vershort}

\begin{verlong}
We shall use the notation introduced in Definition \ref{def.addcol}. We know
that $B$ is an $n\times\left(  n-1\right)  $-matrix (since $B\in
\mathbb{K}^{n\times\left(  n-1\right)  }$); hence, $B_{\bullet,q}$ is an
$n\times1$-matrix. Thus, $\left(  B\mid B_{\bullet,q}\right)  $ is an
$n\times\left(  \left(  n-1\right)  +1\right)  $-matrix. In other words,
$\left(  B\mid B_{\bullet,q}\right)  $ is an $n\times n$-matrix (since
$\left(  n-1\right)  +1=n$). We have%
\begin{equation}
\det\left(  B\mid B_{\bullet,q}\right)  =0
\label{pf.prop.pluecker.pre.row.det=0}%
\end{equation}
\footnote{\textit{Proof of (\ref{pf.prop.pluecker.pre.row.det=0}):} Let
$A=\left(  B\mid B_{\bullet,q}\right)  $. Thus, $A$ is an $n\times n$-matrix
(since $\left(  B\mid B_{\bullet,q}\right)  $ is an $n\times n$-matrix). The
definition of $A_{\bullet,n}$ yields that $A_{\bullet,n}$ is the $n$-th column
of the matrix $A$. Thus,%
\[
A_{\bullet,n}=\left(  \text{the }n\text{-th column of the matrix }A\right)  .
\]
Hence,%
\begin{align}
&  \left(  \text{the }n\text{-th column of the matrix }A\right) \nonumber\\
&  =A_{\bullet,n}=\left(  B\mid B_{\bullet,q}\right)  _{\bullet,n}%
\ \ \ \ \ \ \ \ \ \ \left(  \text{since }A=\left(  B\mid B_{\bullet,q}\right)
\right) \nonumber\\
&  =\left(  B\mid B_{\bullet,q}\right)  _{\bullet,\left(  \left(  n-1\right)
+1\right)  }\ \ \ \ \ \ \ \ \ \ \left(  \text{since }n=\left(  n-1\right)
+1\right) \nonumber\\
&  =B_{\bullet,q}\ \ \ \ \ \ \ \ \ \ \left(
\begin{array}
[c]{c}%
\text{by Proposition \ref{prop.addcol.props1} \textbf{(b)}, applied to
}n-1\text{, }B\text{ and }B_{\bullet,q}\\
\text{instead of }m\text{, }A\text{ and }v
\end{array}
\right)  . \label{pf.prop.pluecker.pre.row.det=0.pf.1}%
\end{align}
\par
On the other hand, $n\in\left\{  1,2,\ldots,n\right\}  $ (since $n$ is a
positive integer) and $q\in\left\{  1,2,\ldots,n-1\right\}  \subseteq\left\{
1,2,\ldots,n\right\}  $. Furthermore, $q\in\left\{  1,2,\ldots,n-1\right\}  $,
thus $q\leq n-1<n$. Hence, $q\neq n$. Hence, $q$ and $n$ are two distinct
elements of $\left\{  1,2,\ldots,n\right\}  $.
\par
But $q\in\left\{  1,2,\ldots,n\right\}  $. Hence, the definition of
$A_{\bullet,q}$ yields that $A_{\bullet,q}$ is the $q$-th column of the matrix
$A$. Thus,%
\[
A_{\bullet,q}=\left(  \text{the }q\text{-th column of the matrix }A\right)  .
\]
Hence,%
\begin{align*}
&  \left(  \text{the }q\text{-th column of the matrix }A\right) \\
&  =A_{\bullet,q}=\left(  B\mid B_{\bullet,q}\right)  _{\bullet,q}%
\ \ \ \ \ \ \ \ \ \ \left(  \text{since }A=\left(  B\mid B_{\bullet,q}\right)
\right) \\
&  =B_{\bullet,q}\ \ \ \ \ \ \ \ \ \ \left(
\begin{array}
[c]{c}%
\text{by Proposition \ref{prop.addcol.props1} \textbf{(a)}, applied to
}n-1\text{, }B\text{ and }B_{\bullet,q}\\
\text{instead of }m\text{, }A\text{ and }v
\end{array}
\right) \\
&  =\left(  \text{the }n\text{-th column of the matrix }A\right)
\ \ \ \ \ \ \ \ \ \ \left(  \text{by
(\ref{pf.prop.pluecker.pre.row.det=0.pf.1})}\right)  .
\end{align*}
Hence, there exist two distinct elements $u$ and $v$ of $\left\{
1,2,\ldots,n\right\}  $ such that
\[
\left(  \text{the }u\text{-th column of the matrix }A\right)  =\left(
\text{the }v\text{-th column of the matrix }A\right)
\]
(namely, $u=q$ and $v=n$). In other words, the matrix $A$ has two equal
columns. Thus, $\det A=0$ (by Exercise \ref{exe.ps4.6} \textbf{(f)}). This
rewrites as $\det\left(  B\mid B_{\bullet,q}\right)  =0$ (since $A=\left(
B\mid B_{\bullet,q}\right)  $). This proves
(\ref{pf.prop.pluecker.pre.row.det=0}).}.
\end{verlong}

But $B_{\bullet,q}$ is the $q$-th column of the matrix $B$ (by the definition
of $B_{\bullet,q}$). Thus,%
\begin{align*}
B_{\bullet,q}  &  =\left(  \text{the }q\text{-th column of the matrix
}B\right) \\
&  =\left(
\begin{array}
[c]{c}%
b_{1,q}\\
b_{2,q}\\
\vdots\\
b_{n,q}%
\end{array}
\right)  \ \ \ \ \ \ \ \ \ \ \left(  \text{since }B=\left(  b_{i,j}\right)
_{1\leq i\leq n,\ 1\leq j\leq n-1}\right) \\
&  =\left(  b_{1,q},b_{2,q},\ldots,b_{n,q}\right)  ^{T}.
\end{align*}
Hence, Proposition \ref{prop.addcol.props2} \textbf{(a)} (applied to $B$,
$B_{\bullet,q}$ and $b_{i,q}$ instead of $A$, $v$ and $v_{i}$) yields%
\[
\det\left(  B\mid B_{\bullet,q}\right)  =\sum_{i=1}^{n}\left(  -1\right)
^{n+i}b_{i,q}\det\left(  B_{\sim i,\bullet}\right)  .
\]
Comparing this with $\det\left(  B\mid B_{\bullet,q}\right)  =0$, we obtain%
\[
0=\sum_{i=1}^{n}\left(  -1\right)  ^{n+i}b_{i,q}\det\left(  B_{\sim i,\bullet
}\right)  .
\]
Multiplying both sides of this equality by $\left(  -1\right)  ^{n}$, we
obtain%
\begin{align*}
0  &  =\left(  -1\right)  ^{n}\sum_{i=1}^{n}\left(  -1\right)  ^{n+i}%
b_{i,q}\det\left(  B_{\sim i,\bullet}\right) \\
&  =\sum_{i=1}^{n}\underbrace{\left(  -1\right)  ^{n}\left(  -1\right)
^{n+i}}_{\substack{=\left(  -1\right)  ^{n+\left(  n+i\right)  }=\left(
-1\right)  ^{i}\\\text{(since }n+\left(  n+i\right)  =2n+i\equiv
i\operatorname{mod}2\text{)}}}\underbrace{b_{i,q}\det\left(  B_{\sim
i,\bullet}\right)  }_{=\det\left(  B_{\sim i,\bullet}\right)  b_{i,q}}\\
&  =\sum_{i=1}^{n}\left(  -1\right)  ^{i}\det\left(  B_{\sim i,\bullet
}\right)  b_{i,q}=\sum_{r=1}^{n}\left(  -1\right)  ^{r}\det\left(  B_{\sim
r,\bullet}\right)  b_{r,q}%
\end{align*}
(here, we renamed the summation index $i$ as $r$). This proves Proposition
\ref{prop.pluecker.pre.row} \textbf{(b)}.

\textbf{(a)} Write the matrix $B$ in the form $B=\left(  b_{i,j}\right)
_{1\leq i\leq n,\ 1\leq j\leq n-1}$. For every $r\in\left\{  1,2,\ldots
,n\right\}  $, we have%
\begin{align*}
B_{r,\bullet}  &  =\left(  \text{the }r\text{-th row of the matrix }B\right)
\\
&  \ \ \ \ \ \ \ \ \ \ \left(
\begin{array}
[c]{c}%
\text{since }B_{r,\bullet}\text{ is the }r\text{-th row of the matrix }B\\
\text{(by the definition of }B_{r,\bullet}\text{)}%
\end{array}
\right) \\
&  =\left(  b_{r,1},b_{r,2},\ldots,b_{r,n-1}\right)
\ \ \ \ \ \ \ \ \ \ \left(  \text{since }B=\left(  b_{i,j}\right)  _{1\leq
i\leq n,\ 1\leq j\leq n-1}\right) \\
&  =\left(  b_{r,j}\right)  _{1\leq i\leq1,\ 1\leq j\leq n-1}.
\end{align*}
Thus,%
\begin{align*}
&  \sum_{r=1}^{n}\left(  -1\right)  ^{r}\det\left(  B_{\sim r,\bullet}\right)
\underbrace{B_{r,\bullet}}_{=\left(  b_{r,j}\right)  _{1\leq i\leq1,\ 1\leq
j\leq n-1}}\\
&  =\sum_{r=1}^{n}\underbrace{\left(  -1\right)  ^{r}\det\left(  B_{\sim
r,\bullet}\right)  \left(  b_{r,j}\right)  _{1\leq i\leq1,\ 1\leq j\leq n-1}%
}_{=\left(  \left(  -1\right)  ^{r}\det\left(  B_{\sim r,\bullet}\right)
b_{r,j}\right)  _{1\leq i\leq1,\ 1\leq j\leq n-1}}=\sum_{r=1}^{n}\left(
\left(  -1\right)  ^{r}\det\left(  B_{\sim r,\bullet}\right)  b_{r,j}\right)
_{1\leq i\leq1,\ 1\leq j\leq n-1}\\
&  =\left(  \underbrace{\sum_{r=1}^{n}\left(  -1\right)  ^{r}\det\left(
B_{\sim r,\bullet}\right)  b_{r,j}}_{\substack{=0\\\text{(by Proposition
\ref{prop.pluecker.pre.row} \textbf{(b)}, applied to }q=j\text{)}}}\right)
_{1\leq i\leq1,\ 1\leq j\leq n-1}=\left(  0\right)  _{1\leq i\leq1,\ 1\leq
j\leq n-1}=0_{1\times\left(  n-1\right)  }.
\end{align*}
This proves Proposition \ref{prop.pluecker.pre.row} \textbf{(a)}.
\end{proof}

Let us next state a variant of Proposition \ref{prop.pluecker.pre.row} where
rows are replaced by columns (and $n$ is renamed as $n+1$):

\begin{proposition}
\label{prop.pluecker.pre.col}Let $n\in\mathbb{N}$. Let $B\in\mathbb{K}%
^{n\times\left(  n+1\right)  }$. Then:

\textbf{(a)} We have%
\[
\sum_{r=1}^{n+1}\left(  -1\right)  ^{r}\det\left(  B_{\bullet,\sim r}\right)
B_{\bullet,r}=0_{n\times1}.
\]

\textbf{(b)} Write the matrix $B$ in the form $B=\left(  b_{i,j}\right)
_{1\leq i\leq n,\ 1\leq j\leq n+1}$. Then,
\[
\sum_{r=1}^{n+1}\left(  -1\right)  ^{r}\det\left(  B_{\bullet,\sim r}\right)
b_{q,r}=0
\]
for every $q\in\left\{  1,2,\ldots,n\right\}  $.
\end{proposition}

\begin{example}
\label{exam.prop.pluecker.pre.col}For this example, set $n=2$ and $B=\left(
\begin{array}
[c]{ccc}%
a & a^{\prime} & a^{\prime\prime}\\
b & b^{\prime} & b^{\prime\prime}%
\end{array}
\right)  $. Then, Proposition \ref{prop.pluecker.pre.col} \textbf{(a)} says
that%
\begin{align*}
&  -\det\left(
\begin{array}
[c]{cc}%
a^{\prime} & a^{\prime\prime}\\
b^{\prime} & b^{\prime\prime}%
\end{array}
\right)  \cdot\left(
\begin{array}
[c]{c}%
a\\
b
\end{array}
\right)  +\det\left(
\begin{array}
[c]{cc}%
a & a^{\prime\prime}\\
b & b^{\prime\prime}%
\end{array}
\right)  \cdot\left(
\begin{array}
[c]{c}%
a^{\prime}\\
b^{\prime}%
\end{array}
\right)  -\det\left(
\begin{array}
[c]{cc}%
a & a^{\prime}\\
b & b^{\prime}%
\end{array}
\right)  \cdot\left(
\begin{array}
[c]{c}%
a^{\prime\prime}\\
b^{\prime\prime}%
\end{array}
\right) \\
&  =0_{2\times1}.
\end{align*}
Proposition \ref{prop.pluecker.pre.col} \textbf{(b)} (applied to $q=2$) yields%
\[
-\det\left(
\begin{array}
[c]{cc}%
a^{\prime} & a^{\prime\prime}\\
b^{\prime} & b^{\prime\prime}%
\end{array}
\right)  \cdot b+\det\left(
\begin{array}
[c]{cc}%
a & a^{\prime\prime}\\
b & b^{\prime\prime}%
\end{array}
\right)  \cdot b^{\prime}-\det\left(
\begin{array}
[c]{cc}%
a & a^{\prime}\\
b & b^{\prime}%
\end{array}
\right)  \cdot b^{\prime\prime}=0.
\]

\end{example}

We shall obtain Proposition \ref{prop.pluecker.pre.col} by applying
Proposition \ref{prop.pluecker.pre.row} to $n+1$ and $B^{T}$
\ \ \ \ \footnote{Recall that $B^{T}$ denotes the transpose of the matrix $B$
(see Definition \ref{def.transpose}).} instead of $n$ and $B$. For this, we
shall need a really simple lemma:

\begin{lemma}
\label{lem.unrows.transpose.1}Let $n\in\mathbb{N}$ and $m\in\mathbb{N}$. Let
$r\in\left\{  1,2,\ldots,m\right\}  $. Let $B\in\mathbb{K}^{n\times m}$. Then,
$\left(  B^{T}\right)  _{\sim r,\bullet}=\left(  B_{\bullet,\sim r}\right)
^{T}$.
\end{lemma}

\begin{vershort}
\begin{proof}
[Proof of Lemma \ref{lem.unrows.transpose.1}.]Taking the transpose of a matrix
turns all its columns into rows. Thus, if we remove the $r$-th \textbf{column}
from $B$ and then take the transpose of the resulting matrix, then we obtain
the same matrix as if we first take the transpose of $B$ and then remove the
$r$-th \textbf{row} from it. Translating this statement into formulas, we
obtain precisely $\left(  B_{\bullet,\sim r}\right)  ^{T}=\left(
B^{T}\right)  _{\sim r,\bullet}$. Thus, Lemma \ref{lem.unrows.transpose.1} is
proven.\footnote{A more formal proof of this could be given using Proposition
\ref{prop.submatrix.easy} \textbf{(e)}.}
\end{proof}
\end{vershort}

\begin{verlong}
\begin{proof}
[Proof of Lemma \ref{lem.unrows.transpose.1}.]We shall use the notations
defined in Definition \ref{def.submatrix} and in Definition \ref{def.rowscols}%
. The definition of $B_{\bullet,\sim r}$ yields
\[
B_{\bullet,\sim r}=\operatorname*{cols}\nolimits_{1,2,\ldots,\widehat{r}%
,\ldots,m}B=\operatorname*{sub}\nolimits_{1,2,\ldots,n}^{1,2,\ldots
,\widehat{r},\ldots,m}B
\]
(by Proposition \ref{prop.submatrix.easy} \textbf{(c)}, applied to $A=B$ and
$v=m-1$ and $\left(  j_{1},j_{2},\ldots,j_{v}\right)  =\left(  1,2,\ldots
,\widehat{r},\ldots,m\right)  $). Hence,%
\begin{equation}
\left(  \underbrace{B_{\bullet,\sim r}}_{=\operatorname*{sub}%
\nolimits_{1,2,\ldots,n}^{1,2,\ldots,\widehat{r},\ldots,m}B}\right)
^{T}=\left(  \operatorname*{sub}\nolimits_{1,2,\ldots,n}^{1,2,\ldots
,\widehat{r},\ldots,m}B\right)  ^{T}=\operatorname*{sub}\nolimits_{1,2,\ldots
,\widehat{r},\ldots,m}^{1,2,\ldots,n}\left(  B^{T}\right)
\label{pf.lem.unrows.transpose.1.1}%
\end{equation}
(by Proposition \ref{prop.submatrix.easy} \textbf{(e)}, applied to $A=B$,
$u=n$, $v=m-1$, $\left(  i_{1},i_{2},\ldots,i_{u}\right)  =\left(
1,2,\ldots,n\right)  $ and $\left(  j_{1},j_{2},\ldots,j_{v}\right)  =\left(
1,2,\ldots,\widehat{r},\ldots,m\right)  $).

On the other hand, $B^{T}$ is an $m\times n$-matrix (since $B$ is an $n\times
m$-matrix (since $B\in\mathbb{K}^{n\times m}$)). Hence, the definition of
$\left(  B^{T}\right)  _{\sim r,\bullet}$ yields%
\[
\left(  B^{T}\right)  _{\sim r,\bullet}=\operatorname*{rows}%
\nolimits_{1,2,\ldots,\widehat{r},\ldots,m}\left(  B^{T}\right)
=\operatorname*{sub}\nolimits_{1,2,\ldots,\widehat{r},\ldots,m}^{1,2,\ldots
,n}\left(  B^{T}\right)
\]
(by Proposition \ref{prop.submatrix.easy} \textbf{(b)}, applied to $m$, $n$,
$B^{T}$, $m-1$ and $\left(  1,2,\ldots,\widehat{r},\ldots,m\right)  $ instead
of $n$, $m$, $A$, $u$ and $\left(  i_{1},i_{2},\ldots,i_{u}\right)  $).
Comparing this with (\ref{pf.lem.unrows.transpose.1.1}), we obtain $\left(
B^{T}\right)  _{\sim r,\bullet}=\left(  B_{\bullet,\sim r}\right)  ^{T}$.
Thus, Lemma \ref{lem.unrows.transpose.1} is proven.
\end{proof}
\end{verlong}

\begin{proof}
[Proof of Proposition \ref{prop.pluecker.pre.col}.]\textbf{(b)} Let
$q\in\left\{  1,2,\ldots,n\right\}  $. Thus, \newline$q\in\left\{
1,2,\ldots,n\right\}  =\left\{  1,2,\ldots,\left(  n+1\right)  -1\right\}  $
(since $n=\left(  n+1\right)  -1$).

We have $B=\left(  b_{i,j}\right)  _{1\leq i\leq n,\ 1\leq j\leq n+1}$. Thus,
the definition of $B^{T}$ yields%
\[
B^{T}=\left(  b_{j,i}\right)  _{1\leq i\leq n+1,\ 1\leq j\leq n}=\left(
b_{j,i}\right)  _{1\leq i\leq n+1,\ 1\leq j\leq\left(  n+1\right)  -1}%
\]
(since $n=\left(  n+1\right)  -1$). Also, $B^{T}=\left(  b_{j,i}\right)
_{1\leq i\leq n+1,\ 1\leq j\leq\left(  n+1\right)  -1}\in\mathbb{K}^{\left(
n+1\right)  \times\left(  \left(  n+1\right)  -1\right)  }$. Thus, Proposition
\ref{prop.pluecker.pre.row} \textbf{(b)} (applied to $n+1$, $B^{T}$ and
$b_{j,i}$ instead of $n$, $B$ and $b_{i,j}$) yields%
\begin{equation}
\sum_{r=1}^{n+1}\left(  -1\right)  ^{r}\det\left(  \left(  B^{T}\right)
_{\sim r,\bullet}\right)  b_{q,r}=0. \label{pf.prop.pluecker.pre.col.b.1}%
\end{equation}
But every $r\in\left\{  1,2,\ldots,n+1\right\}  $ satisfies%
\begin{equation}
\det\left(  \left(  B^{T}\right)  _{\sim r,\bullet}\right)  =\det\left(
B_{\bullet,\sim r}\right)  \label{pf.prop.pluecker.pre.col.b.2}%
\end{equation}
\footnote{\textit{Proof of (\ref{pf.prop.pluecker.pre.col.b.2}):} Let
$r\in\left\{  1,2,\ldots,n+1\right\}  $. Then, $B_{\bullet,\sim r}%
\in\mathbb{K}^{n\times n}$ (since $B\in\mathbb{K}^{n\times\left(  n+1\right)
}$). In other words, $B_{\bullet,\sim r}$ is an $n\times n$-matrix. Thus,
Exercise \ref{exe.ps4.4} (applied to $A=B_{\bullet,\sim r}$) yields
$\det\left(  \left(  B_{\bullet,\sim r}\right)  ^{T}\right)  =\det\left(
B_{\bullet,\sim r}\right)  $.
\par
But Lemma \ref{lem.unrows.transpose.1} (applied to $m=n+1$) yields $\left(
B^{T}\right)  _{\sim r,\bullet}=\left(  B_{\bullet,\sim r}\right)  ^{T}$.
Thus, $\det\left(  \underbrace{\left(  B^{T}\right)  _{\sim r,\bullet}%
}_{=\left(  B_{\bullet,\sim r}\right)  ^{T}}\right)  =\det\left(  \left(
B_{\bullet,\sim r}\right)  ^{T}\right)  =\det\left(  B_{\bullet,\sim
r}\right)  $. This proves (\ref{pf.prop.pluecker.pre.col.b.2}).}. Hence,
(\ref{pf.prop.pluecker.pre.col.b.1}) yields%
\[
0=\sum_{r=1}^{n+1}\left(  -1\right)  ^{r}\underbrace{\det\left(  \left(
B^{T}\right)  _{\sim r,\bullet}\right)  }_{\substack{=\det\left(
B_{\bullet,\sim r}\right)  \\\text{(by (\ref{pf.prop.pluecker.pre.col.b.2}))}%
}}b_{q,r}=\sum_{r=1}^{n+1}\left(  -1\right)  ^{r}\det\left(  B_{\bullet,\sim
r}\right)  b_{q,r}.
\]
This proves Proposition \ref{prop.pluecker.pre.col} \textbf{(b)}.

\textbf{(a)} Write the matrix $B$ in the form $B=\left(  b_{i,j}\right)
_{1\leq i\leq n,\ 1\leq j\leq n+1}$. For every $r\in\left\{  1,2,\ldots
,n+1\right\}  $, we have%
\begin{align*}
B_{\bullet,r}  &  =\left(  \text{the }r\text{-th column of the matrix
}B\right) \\
&  \ \ \ \ \ \ \ \ \ \ \left(
\begin{array}
[c]{c}%
\text{since }B_{\bullet,r}\text{ is the }r\text{-th column of the matrix }B\\
\text{(by the definition of }B_{\bullet,r}\text{)}%
\end{array}
\right) \\
&  =\left(
\begin{array}
[c]{c}%
b_{1,r}\\
b_{2,r}\\
\vdots\\
b_{n,r}%
\end{array}
\right)  \ \ \ \ \ \ \ \ \ \ \left(  \text{since }B=\left(  b_{i,j}\right)
_{1\leq i\leq n,\ 1\leq j\leq n+1}\right) \\
&  =\left(  b_{i,r}\right)  _{1\leq i\leq n,\ 1\leq j\leq1}.
\end{align*}
Thus,%
\begin{align*}
&  \sum_{r=1}^{n+1}\left(  -1\right)  ^{r}\det\left(  B_{\bullet,\sim
r}\right)  \underbrace{B_{\bullet,r}}_{=\left(  b_{i,r}\right)  _{1\leq i\leq
n,\ 1\leq j\leq1}}\\
&  =\sum_{r=1}^{n+1}\underbrace{\left(  -1\right)  ^{r}\det\left(
B_{\bullet,\sim r}\right)  \left(  b_{i,r}\right)  _{1\leq i\leq n,\ 1\leq
j\leq1}}_{=\left(  \left(  -1\right)  ^{r}\det\left(  B_{\bullet,\sim
r}\right)  b_{i,r}\right)  _{1\leq i\leq n,\ 1\leq j\leq1}}=\sum_{r=1}%
^{n+1}\left(  \left(  -1\right)  ^{r}\det\left(  B_{\bullet,\sim r}\right)
b_{i,r}\right)  _{1\leq i\leq n,\ 1\leq j\leq1}\\
&  =\left(  \underbrace{\sum_{r=1}^{n+1}\left(  -1\right)  ^{r}\det\left(
B_{\bullet,\sim r}\right)  b_{i,r}}_{\substack{=0\\\text{(by Proposition
\ref{prop.pluecker.pre.col} \textbf{(b)}, applied to }q=i\text{)}}}\right)
_{1\leq i\leq n,\ 1\leq j\leq1}=\left(  0\right)  _{1\leq i\leq n,\ 1\leq
j\leq1}=0_{n\times1}.
\end{align*}
This proves Proposition \ref{prop.pluecker.pre.col} \textbf{(a)}.
\end{proof}

\begin{remark}
\label{rmk.pluecker.pre.col.cramer}Proposition \ref{prop.pluecker.pre.col}
\textbf{(a)} can be viewed as a restatement of Cramer's rule (Theorem
\ref{thm.cramer} \textbf{(a)}). More precisely, it is easy to derive one of
these two facts from the other (although neither of the two is difficult to
prove to begin with). Let us sketch one direction of this argument: namely,
let us derive Theorem \ref{thm.cramer} \textbf{(a)} from Proposition
\ref{prop.pluecker.pre.col} \textbf{(a)}.

Indeed, let $n$, $A$, $b=\left(  b_{1},b_{2},\ldots,b_{n}\right)  ^{T}$ and
$A_{j}^{\#}$ be as in Theorem \ref{thm.cramer} \textbf{(a)}. We want to prove
that $A\cdot\left(  \det\left(  A_{1}^{\#}\right)  ,\det\left(  A_{2}%
^{\#}\right)  ,\ldots,\det\left(  A_{n}^{\#}\right)  \right)  ^{T}=\det A\cdot
b$.

Let $B=\left(  A\mid b\right)  $ (using the notations of Definition
\ref{def.addcol}); this is an $n\times\left(  n+1\right)  $-matrix.

Fix $r\in\left\{  1,2,\ldots,n\right\}  $. The matrix $B_{\bullet,\sim r}$
differs from the matrix $A_{r}^{\#}$ only in the order of its columns: More precisely,

\begin{itemize}
\item the matrix $B_{\bullet,\sim r}$ is obtained from the matrix $A$ by
removing the $r$-th column and attaching the column vector $b$ to the right
edge, whereas

\item the matrix $A_{r}^{\#}$ is obtained from the matrix $A$ by replacing the
$r$-th column by the column vector $b$.
\end{itemize}

Thus, the matrix $B_{\bullet,\sim r}$ can be obtained from the matrix
$A_{r}^{\#}$ by first swapping the $r$-th and $\left(  r+1\right)  $-th
columns, then swapping the $\left(  r+1\right)  $-th and $\left(  r+2\right)
$-th columns, etc., until finally swapping the $\left(  n-1\right)  $-th and
$n$-th columns. Each of these swaps multiplies the determinant by $-1$ (by
Exercise \ref{exe.ps4.6} \textbf{(b)}); thus, our sequence of swaps multiplies
the determinant by $\left(  -1\right)  ^{n-r}=\left(  -1\right)  ^{n+r}$.
Hence,
\begin{equation}
\det\left(  B_{\bullet,\sim r}\right)  =\left(  -1\right)  ^{n+r}\det\left(
A_{r}^{\#}\right)  . \label{eq.rmk.pluecker.pre.col.cramer.1}%
\end{equation}

Now, forget that we fixed $r$. It is easy to see that every $\left(
v_{1},v_{2},\ldots,v_{n}\right)  ^{T}\in\mathbb{K}^{1\times n}$ satisfies%
\[
A\cdot\left(  v_{1},v_{2},\ldots,v_{n}\right)  ^{T}=\sum_{r=1}^{n}%
v_{r}A_{\bullet,r}.
\]
Applying this to $\left(  v_{1},v_{2},\ldots,v_{n}\right)  ^{T}=\left(
\det\left(  A_{1}^{\#}\right)  ,\det\left(  A_{2}^{\#}\right)  ,\ldots
,\det\left(  A_{n}^{\#}\right)  \right)  ^{T}$, we obtain%
\begin{equation}
A\cdot\left(  \det\left(  A_{1}^{\#}\right)  ,\det\left(  A_{2}^{\#}\right)
,\ldots,\det\left(  A_{n}^{\#}\right)  \right)  ^{T}=\sum_{r=1}^{n}\det\left(
A_{r}^{\#}\right)  A_{\bullet,r}. \label{eq.rmk.pluecker.pre.col.cramer.2}%
\end{equation}
But Proposition \ref{prop.pluecker.pre.col} \textbf{(a)} yields%
\begin{align*}
0_{n\times1}  &  =\sum_{r=1}^{n+1}\left(  -1\right)  ^{r}\det\left(
B_{\bullet,\sim r}\right)  B_{\bullet,r}\\
&  =\sum_{r=1}^{n}\left(  -1\right)  ^{r}\underbrace{\det\left(
B_{\bullet,\sim r}\right)  }_{\substack{=\left(  -1\right)  ^{n+r}\det\left(
A_{r}^{\#}\right)  \\\text{(by (\ref{eq.rmk.pluecker.pre.col.cramer.1}))}%
}}\underbrace{B_{\bullet,r}}_{\substack{=A_{\bullet,r}\\\text{(since
}B=\left(  A\mid b\right)  \text{ and }r\leq n\text{)}}}\\
&  \ \ \ \ \ \ \ \ \ \ +\left(  -1\right)  ^{n+1}\det\left(
\underbrace{B_{\bullet,\sim\left(  n+1\right)  }}_{\substack{=A\\\text{(since
}B=\left(  A\mid b\right)  \text{)}}}\right)  \underbrace{B_{\bullet,n+1}%
}_{\substack{=b\\\text{(since }B=\left(  A\mid b\right)  \text{)}}}\\
&  =\sum_{r=1}^{n}\underbrace{\left(  -1\right)  ^{r}\left(  -1\right)
^{n+r}}_{=\left(  -1\right)  ^{n}}\det\left(  A_{r}^{\#}\right)  A_{\bullet
,r}+\underbrace{\left(  -1\right)  ^{n+1}}_{=-\left(  -1\right)  ^{n}}\det
A\cdot b\\
&  =\sum_{r=1}^{n}\left(  -1\right)  ^{n}\det\left(  A_{r}^{\#}\right)
A_{\bullet,r}-\left(  -1\right)  ^{n}\det A\cdot b\\
&  =\left(  -1\right)  ^{n}\left(  \sum_{r=1}^{n}\det\left(  A_{r}%
^{\#}\right)  A_{\bullet,r}-\det A\cdot b\right)  .
\end{align*}
Multiplying both sides of this equality by $\left(  -1\right)  ^{n}$, we
obtain%
\begin{align*}
0_{n\times1}  &  =\underbrace{\left(  -1\right)  ^{n}\left(  -1\right)  ^{n}%
}_{=1}\left(  \sum_{r=1}^{n}\det\left(  A_{r}^{\#}\right)  A_{\bullet,r}-\det
A\cdot b\right) \\
&  =\sum_{r=1}^{n}\det\left(  A_{r}^{\#}\right)  A_{\bullet,r}-\det A\cdot b.
\end{align*}
Hence,
\[
\det A\cdot b=\sum_{r=1}^{n}\det\left(  A_{r}^{\#}\right)  A_{\bullet
,r}=A\cdot\left(  \det\left(  A_{1}^{\#}\right)  ,\det\left(  A_{2}%
^{\#}\right)  ,\ldots,\det\left(  A_{n}^{\#}\right)  \right)  ^{T}%
\]
(by (\ref{eq.rmk.pluecker.pre.col.cramer.2})). Thus, we have derived Theorem
\ref{thm.cramer} \textbf{(a)} from Proposition \ref{prop.pluecker.pre.col}
\textbf{(a)}. Essentially the same argument (but read backwards) can be used
to derive Proposition \ref{prop.pluecker.pre.col} \textbf{(a)} from Theorem
\ref{thm.cramer} \textbf{(a)}.
\end{remark}

Now, we can easily prove the \textit{Pl\"{u}cker identity}:

\begin{theorem}
\label{thm.pluecker.plu}Let $n$ be a positive integer. Let $A\in
\mathbb{K}^{n\times\left(  n-1\right)  }$ and $B\in\mathbb{K}^{n\times\left(
n+1\right)  }$. Then,%
\[
\sum_{r=1}^{n+1}\left(  -1\right)  ^{r}\det\left(  A\mid B_{\bullet,r}\right)
\det\left(  B_{\bullet,\sim r}\right)  =0
\]
(where we are using the notations from Definition \ref{def.addcol} and from
Definition \ref{def.unrows}).
\end{theorem}

\begin{example}
\label{exam.thm.pluecker.plu}If $n=3$, $A=\left(
\begin{array}
[c]{cc}%
a & a^{\prime}\\
b & b^{\prime}\\
c & c^{\prime}%
\end{array}
\right)  $ and $B=\left(
\begin{array}
[c]{cccc}%
x_{1} & x_{2} & x_{3} & x_{4}\\
y_{1} & y_{2} & y_{3} & y_{4}\\
z_{1} & z_{2} & z_{3} & z_{4}%
\end{array}
\right)  $, then Theorem \ref{thm.pluecker.plu} says that%
\begin{align*}
&  -\det\left(
\begin{array}
[c]{ccc}%
a & a^{\prime} & x_{1}\\
b & b^{\prime} & y_{1}\\
c & c^{\prime} & z_{1}%
\end{array}
\right)  \det\left(
\begin{array}
[c]{ccc}%
x_{2} & x_{3} & x_{4}\\
y_{2} & y_{3} & y_{4}\\
z_{2} & z_{3} & z_{4}%
\end{array}
\right) \\
&  \ \ \ \ \ \ \ \ \ \ +\det\left(
\begin{array}
[c]{ccc}%
a & a^{\prime} & x_{2}\\
b & b^{\prime} & y_{2}\\
c & c^{\prime} & z_{2}%
\end{array}
\right)  \det\left(
\begin{array}
[c]{ccc}%
x_{1} & x_{3} & x_{4}\\
y_{1} & y_{3} & y_{4}\\
z_{1} & z_{3} & z_{4}%
\end{array}
\right) \\
&  \ \ \ \ \ \ \ \ \ \ -\det\left(
\begin{array}
[c]{ccc}%
a & a^{\prime} & x_{3}\\
b & b^{\prime} & y_{3}\\
c & c^{\prime} & z_{3}%
\end{array}
\right)  \det\left(
\begin{array}
[c]{ccc}%
x_{1} & x_{2} & x_{4}\\
y_{1} & y_{2} & y_{4}\\
z_{1} & z_{2} & z_{4}%
\end{array}
\right) \\
&  \ \ \ \ \ \ \ \ \ \ +\det\left(
\begin{array}
[c]{ccc}%
a & a^{\prime} & x_{4}\\
b & b^{\prime} & y_{4}\\
c & c^{\prime} & z_{4}%
\end{array}
\right)  \det\left(
\begin{array}
[c]{ccc}%
x_{1} & x_{2} & x_{3}\\
y_{1} & y_{2} & y_{3}\\
z_{1} & z_{2} & z_{3}%
\end{array}
\right) \\
&  =0.
\end{align*}

\end{example}

\begin{proof}
[Proof of Theorem \ref{thm.pluecker.plu}.]Write the matrix $B$ in the form
$B=\left(  b_{i,j}\right)  _{1\leq i\leq n,\ 1\leq j\leq n+1}$.

Let $r\in\left\{  1,2,\ldots,n+1\right\}  $. Then,%
\begin{align*}
B_{\bullet,r}  &  =\left(  \text{the }r\text{-th column of the matrix
}B\right) \\
&  \ \ \ \ \ \ \ \ \ \ \left(
\begin{array}
[c]{c}%
\text{since }B_{\bullet,r}\text{ is the }r\text{-th column of the matrix }B\\
\text{(by the definition of }B_{\bullet,r}\text{)}%
\end{array}
\right) \\
&  =\left(
\begin{array}
[c]{c}%
b_{1,r}\\
b_{2,r}\\
\vdots\\
b_{n,r}%
\end{array}
\right)  \ \ \ \ \ \ \ \ \ \ \left(  \text{since }B=\left(  b_{i,j}\right)
_{1\leq i\leq n,\ 1\leq j\leq n+1}\right) \\
&  =\left(  b_{1,r},b_{2,r},\ldots,b_{n,r}\right)  ^{T}.
\end{align*}
Hence, Proposition \ref{prop.addcol.props2} \textbf{(a)} (applied to
$B_{\bullet,r}$ and $b_{i,r}$ instead of $v$ and $v_{i}$) shows that%
\begin{equation}
\det\left(  A\mid B_{\bullet,r}\right)  =\sum_{i=1}^{n}\left(  -1\right)
^{n+i}b_{i,r}\det\left(  A_{\sim i,\bullet}\right)  .
\label{pf.thm.pluecker.plu.1}%
\end{equation}

Now, forget that we fixed $r$. We thus have proven
(\ref{pf.thm.pluecker.plu.1}) for each $r\in\left\{  1,2,\ldots,n+1\right\}
$. Now,%
\begin{align*}
&  \sum_{r=1}^{n+1}\left(  -1\right)  ^{r}\underbrace{\det\left(  A\mid
B_{\bullet,r}\right)  }_{\substack{=\sum_{i=1}^{n}\left(  -1\right)
^{n+i}b_{i,r}\det\left(  A_{\sim i,\bullet}\right)  \\\text{(by
(\ref{pf.thm.pluecker.plu.1}))}}}\det\left(  B_{\bullet,\sim r}\right) \\
&  =\sum_{r=1}^{n+1}\left(  -1\right)  ^{r}\left(  \sum_{i=1}^{n}\left(
-1\right)  ^{n+i}b_{i,r}\det\left(  A_{\sim i,\bullet}\right)  \right)
\det\left(  B_{\bullet,\sim r}\right) \\
&  =\underbrace{\sum_{r=1}^{n+1}\sum_{i=1}^{n}}_{=\sum_{i=1}^{n}\sum
_{r=1}^{n+1}}\underbrace{\left(  -1\right)  ^{r}\left(  -1\right)
^{n+i}b_{i,r}\det\left(  A_{\sim i,\bullet}\right)  \det\left(  B_{\bullet
,\sim r}\right)  }_{=\left(  -1\right)  ^{r}\det\left(  B_{\bullet,\sim
r}\right)  b_{i,r}\cdot\left(  -1\right)  ^{n+i}\det\left(  A_{\sim i,\bullet
}\right)  }\\
&  =\sum_{i=1}^{n}\sum_{r=1}^{n+1}\left(  -1\right)  ^{r}\det\left(
B_{\bullet,\sim r}\right)  b_{i,r}\cdot\left(  -1\right)  ^{n+i}\det\left(
A_{\sim i,\bullet}\right) \\
&  =\sum_{i=1}^{n}\underbrace{\left(  \sum_{r=1}^{n+1}\left(  -1\right)
^{r}\det\left(  B_{\bullet,\sim r}\right)  b_{i,r}\right)  }%
_{\substack{=0\\\text{(by Proposition \ref{prop.pluecker.pre.col}
\textbf{(b)},}\\\text{applied to }q=i\text{)}}}\left(  -1\right)  ^{n+i}%
\det\left(  A_{\sim i,\bullet}\right) \\
&  =\sum_{i=1}^{n}0\left(  -1\right)  ^{n+i}\det\left(  A_{\sim i,\bullet
}\right)  =0.
\end{align*}
This proves Theorem \ref{thm.pluecker.plu}.
\end{proof}

\begin{remark}
Theorem \ref{thm.pluecker.plu} (at least in the case when $\mathbb{K}$ is a
field) is essentially equivalent to \cite[\S 9.1, Exercise 1]{Fulton-Young},
to \cite[(QR)]{KleLak72}, to \cite[Theorem 3.4.11 (the \textquotedblleft
necessary\textquotedblright\ part)]{Jacobs10}, and to \cite[Proposition 3.3.2
(the \textquotedblleft only if\textquotedblright\ part)]{Lampe}.
\end{remark}

\begin{verlong}
\begin{remark}
\label{rmk.thm.pluecker.plu.GR}Theorem \ref{thm.pluecker.plu} is equivalent to
Proposition 14.4 in the detailed version of the paper\newline[GriRob14] Darij
Grinberg, Tom Roby, \textit{\textit{Iterative properties of birational
rowmotion}}, \href{http://arxiv.org/abs/1402.6178v6}{arXiv:1402.6178v6}.

More precisely, Theorem \ref{thm.pluecker.plu} is exactly Proposition 14.4 in
the detailed version of [GriRob14] (applied to $u=n$ and $w_{i}=B_{\bullet,i}%
$). Conversely, Proposition 14.4 in the detailed version of [GriRob14] is
precisely Theorem \ref{thm.pluecker.plu} (applied to $n=u$ and $B=\left(
w_{1}\mid w_{2}\mid\cdots\mid w_{u+1}\right)  $, where we are using the
notations of [GriRob14]).
\end{remark}
\end{verlong}

\begin{exercise}
\label{exe.pluecker.rederive-AC}Use Theorem \ref{thm.pluecker.plu} to give a
new proof of Proposition \ref{prop.desnanot.AC}.
\end{exercise}

\subsection{\label{sect.laplace}Laplace expansion in multiple rows/columns}

In this section, we shall see a (somewhat unwieldy, but classical and
important) generalization of Theorem \ref{thm.laplace.gen}. First, we shall
need some notations:

\begin{definition}
\label{def.sect.laplace.notations}Throughout Section \ref{sect.laplace}, we
shall use the following notations:

\begin{itemize}
\item If $I$ is a finite set of integers, then $\sum I$ shall denote the sum
of all elements of $I$. (Thus, $\sum I=\sum_{i\in I}i$.)

\item If $I$ is a finite set of integers, then $w\left(  I\right)  $ shall
denote the list of all elements of $I$ in increasing order (with no
repetitions). (See Definition \ref{def.ind.inclist} for the formal definition
of this list.) (For example, $w\left(  \left\{  3,4,8\right\}  \right)
=\left(  3,4,8\right)  $.)
\end{itemize}

We shall also use the notation introduced in Definition \ref{def.submatrix}.
If $n$, $m$, $A$, $\left(  i_{1},i_{2},\ldots,i_{u}\right)  $ and $\left(
j_{1},j_{2},\ldots,j_{v}\right)  $ are as in Definition \ref{def.submatrix},
then we shall use the notation $\operatorname*{sub}\nolimits_{\left(
i_{1},i_{2},\ldots,i_{u}\right)  }^{\left(  j_{1},j_{2},\ldots,j_{v}\right)
}A$ as a synonym for $\operatorname*{sub}\nolimits_{i_{1},i_{2},\ldots,i_{u}%
}^{j_{1},j_{2},\ldots,j_{v}}A$.
\end{definition}

A consequence of this definition is that if $A$ is an $n\times m$-matrix, and
if $U$ is a subset of $\left\{  1,2,\ldots,n\right\}  $, and if $V$ is a
subset of $\left\{  1,2,\ldots,m\right\}  $, then $\operatorname*{sub}%
\nolimits_{w\left(  U\right)  }^{w\left(  V\right)  }A$ is a well-defined
$\left\vert U\right\vert \times\left\vert V\right\vert $%
-matrix\footnote{because $w\left(  U\right)  $ is a list of $\left\vert
U\right\vert $ elements of $\left\{  1,2,\ldots,n\right\}  $, and because
$w\left(  V\right)  $ is a list of $\left\vert V\right\vert $ elements of
$\left\{  1,2,\ldots,m\right\}  $} (actually, a submatrix of $A$).

\begin{example}
If $n=3$ and $m=4$ and $A=\left(
\begin{array}
[c]{cccc}%
a & b & c & d\\
a^{\prime} & b^{\prime} & c^{\prime} & d^{\prime}\\
a^{\prime\prime} & b^{\prime\prime} & c^{\prime\prime} & d^{\prime\prime}%
\end{array}
\right)  $, then $\operatorname*{sub}\nolimits_{w\left(  \left\{  2,3\right\}
\right)  }^{w\left(  \left\{  1,3,4\right\}  \right)  }A=\left(
\begin{array}
[c]{ccc}%
a^{\prime} & c^{\prime} & d^{\prime}\\
a^{\prime\prime} & c^{\prime\prime} & d^{\prime\prime}%
\end{array}
\right)  $.
\end{example}

The following fact is obvious from the definition of $w\left(  I\right)  $:

\begin{proposition}
\label{prop.sect.laplace.notations.w(I)}Let $I$ be a finite set of integers.
Then, $w\left(  I\right)  $ is an $\left\vert I\right\vert $-tuple of elements
of $I$.
\end{proposition}

\begin{verlong}
\begin{proof}
[Proof of Proposition \ref{prop.sect.laplace.notations.w(I)}.]The definition
of $w\left(  I\right)  $ says that $w\left(  I\right)  $ is the list of all
elements of $I$ in increasing order (with no repetitions). Thus, $w\left(
I\right)  $ is a list of length $\left\vert I\right\vert $. In other words,
$w\left(  I\right)  $ is an $\left\vert I\right\vert $-tuple. Moreover,
$w\left(  I\right)  $ is a list of elements of $I$\ \ \ \ \footnote{since
$w\left(  I\right)  $ is the list of all elements of $I$ in increasing order
(with no repetitions)}. Hence, $w\left(  I\right)  $ is an $\left\vert
I\right\vert $-tuple of elements of $I$ (since $w\left(  I\right)  $ is an
$\left\vert I\right\vert $-tuple). This proves Proposition
\ref{prop.sect.laplace.notations.w(I)}.
\end{proof}
\end{verlong}

\begin{theorem}
\label{thm.det.laplace-multi}Let $n\in\mathbb{N}$. Let $A\in\mathbb{K}%
^{n\times n}$. For any subset $I$ of $\left\{  1,2,\ldots,n\right\}  $, we let
$\widetilde{I}$ denote the complement $\left\{  1,2,\ldots,n\right\}
\setminus I$ of $I$. (For instance, if $n=4$ and $I=\left\{  1,4\right\}  $,
then $\widetilde{I}=\left\{  2,3\right\}  $.)

\textbf{(a)} For every subset $P$ of $\left\{  1,2,\ldots,n\right\}  $, we
have%
\[
\det A=\sum_{\substack{Q\subseteq\left\{  1,2,\ldots,n\right\}  ;\\\left\vert
Q\right\vert =\left\vert P\right\vert }}\left(  -1\right)  ^{\sum P+\sum
Q}\det\left(  \operatorname*{sub}\nolimits_{w\left(  P\right)  }^{w\left(
Q\right)  }A\right)  \det\left(  \operatorname*{sub}\nolimits_{w\left(
\widetilde{P}\right)  }^{w\left(  \widetilde{Q}\right)  }A\right)  .
\]

\textbf{(b)} For every subset $Q$ of $\left\{  1,2,\ldots,n\right\}  $, we
have%
\[
\det A=\sum_{\substack{P\subseteq\left\{  1,2,\ldots,n\right\}  ;\\\left\vert
P\right\vert =\left\vert Q\right\vert }}\left(  -1\right)  ^{\sum P+\sum
Q}\det\left(  \operatorname*{sub}\nolimits_{w\left(  P\right)  }^{w\left(
Q\right)  }A\right)  \det\left(  \operatorname*{sub}\nolimits_{w\left(
\widetilde{P}\right)  }^{w\left(  \widetilde{Q}\right)  }A\right)  .
\]

\end{theorem}

Theorem \ref{thm.det.laplace-multi} is actually a generalization of Theorem
\ref{thm.laplace.gen}, known as \textquotedblleft Laplace expansion in
multiple rows (resp. columns)\textquotedblright. It appears (for example) in
\cite[Theorem 3.61]{Gill} and in \cite[Theorem 2.4.1]{Prasolov}.\footnote{Of
course, parts \textbf{(a)} and \textbf{(b)} of Theorem
\ref{thm.det.laplace-multi} are easily seen to be equivalent; thus, many
authors confine themselves to only stating one of them. For example, Theorem
\ref{thm.det.laplace-multi} is \cite[Lemma A.1 (f)]{CaSoSp12}.} Theorem
\ref{thm.laplace.gen} \textbf{(a)} can be recovered from Theorem
\ref{thm.det.laplace-multi} \textbf{(a)} by setting $P=\left\{  p\right\}  $;
similarly for the part \textbf{(b)}.

\begin{example}
Let us see what Theorem \ref{thm.det.laplace-multi} \textbf{(a)} says in a
simple case. For this example, set $n=4$ and $A=\left(  a_{i,j}\right)
_{1\leq i\leq4,\ 1\leq j\leq4}$. Also, set $P=\left\{  1,4\right\}
\subseteq\left\{  1,2,3,4\right\}  $; thus, $w\left(  P\right)  =\left(
1,4\right)  $, $\sum P=1+4=5$, $\left\vert P\right\vert =2$, $\widetilde{P}%
=\left\{  2,3\right\}  $ and $w\left(  \widetilde{P}\right)  =\left(
2,3\right)  $. Now, Theorem \ref{thm.det.laplace-multi} \textbf{(a)} says that%
\begin{align*}
&  \det A\\
&  =\sum_{\substack{Q\subseteq\left\{  1,2,3,4\right\}  ;\\\left\vert
Q\right\vert =2}}\left(  -1\right)  ^{5+\sum Q}\det\left(  \operatorname*{sub}%
\nolimits_{1,4}^{w\left(  Q\right)  }A\right)  \det\left(  \operatorname*{sub}%
\nolimits_{2,3}^{w\left(  \widetilde{Q}\right)  }A\right) \\
&  =\left(  -1\right)  ^{5+\left(  1+2\right)  }\det\left(
\operatorname*{sub}\nolimits_{1,4}^{1,2}A\right)  \det\left(
\operatorname*{sub}\nolimits_{2,3}^{3,4}A\right) \\
&  \ \ \ \ \ \ \ \ \ \ +\left(  -1\right)  ^{5+\left(  1+3\right)  }%
\det\left(  \operatorname*{sub}\nolimits_{1,4}^{1,3}A\right)  \det\left(
\operatorname*{sub}\nolimits_{2,3}^{2,4}A\right) \\
&  \ \ \ \ \ \ \ \ \ \ +\left(  -1\right)  ^{5+\left(  1+4\right)  }%
\det\left(  \operatorname*{sub}\nolimits_{1,4}^{1,4}A\right)  \det\left(
\operatorname*{sub}\nolimits_{2,3}^{2,3}A\right) \\
&  \ \ \ \ \ \ \ \ \ \ +\left(  -1\right)  ^{5+\left(  2+3\right)  }%
\det\left(  \operatorname*{sub}\nolimits_{1,4}^{2,3}A\right)  \det\left(
\operatorname*{sub}\nolimits_{2,3}^{1,4}A\right) \\
&  \ \ \ \ \ \ \ \ \ \ +\left(  -1\right)  ^{5+\left(  2+4\right)  }%
\det\left(  \operatorname*{sub}\nolimits_{1,4}^{2,4}A\right)  \det\left(
\operatorname*{sub}\nolimits_{2,3}^{1,3}A\right) \\
&  \ \ \ \ \ \ \ \ \ \ +\left(  -1\right)  ^{5+\left(  3+4\right)  }%
\det\left(  \operatorname*{sub}\nolimits_{1,4}^{3,4}A\right)  \det\left(
\operatorname*{sub}\nolimits_{2,3}^{1,2}A\right) \\
&  =\det\left(
\begin{array}
[c]{cc}%
a_{1,1} & a_{1,2}\\
a_{4,1} & a_{4,2}%
\end{array}
\right)  \det\left(
\begin{array}
[c]{cc}%
a_{2,3} & a_{2,4}\\
a_{3,3} & a_{3,4}%
\end{array}
\right) \\
&  \ \ \ \ \ \ \ \ \ \ -\det\left(
\begin{array}
[c]{cc}%
a_{1,1} & a_{1,3}\\
a_{4,1} & a_{4,3}%
\end{array}
\right)  \det\left(
\begin{array}
[c]{cc}%
a_{2,2} & a_{2,4}\\
a_{3,2} & a_{3,4}%
\end{array}
\right) \\
&  \ \ \ \ \ \ \ \ \ \ +\det\left(
\begin{array}
[c]{cc}%
a_{1,1} & a_{1,4}\\
a_{4,1} & a_{4,4}%
\end{array}
\right)  \det\left(
\begin{array}
[c]{cc}%
a_{2,2} & a_{2,3}\\
a_{3,2} & a_{3,3}%
\end{array}
\right) \\
&  \ \ \ \ \ \ \ \ \ \ +\det\left(
\begin{array}
[c]{cc}%
a_{1,2} & a_{1,3}\\
a_{4,2} & a_{4,3}%
\end{array}
\right)  \det\left(
\begin{array}
[c]{cc}%
a_{2,1} & a_{2,4}\\
a_{3,1} & a_{3,4}%
\end{array}
\right) \\
&  \ \ \ \ \ \ \ \ \ \ -\det\left(
\begin{array}
[c]{cc}%
a_{1,2} & a_{1,4}\\
a_{4,2} & a_{4,4}%
\end{array}
\right)  \det\left(
\begin{array}
[c]{cc}%
a_{2,1} & a_{2,3}\\
a_{3,1} & a_{3,3}%
\end{array}
\right) \\
&  \ \ \ \ \ \ \ \ \ \ +\det\left(
\begin{array}
[c]{cc}%
a_{1,3} & a_{1,4}\\
a_{4,3} & a_{4,4}%
\end{array}
\right)  \det\left(
\begin{array}
[c]{cc}%
a_{2,1} & a_{2,2}\\
a_{3,1} & a_{3,2}%
\end{array}
\right)  .
\end{align*}

\end{example}

The following lemma will play a crucial role in our proof of Theorem
\ref{thm.det.laplace-multi} (similar to the role that Lemma
\ref{lem.laplace.Apq} played in our proof of Theorem \ref{thm.laplace.gen}):

\begin{lemma}
\label{lem.det.laplace-multi.Apq}Let $n\in\mathbb{N}$. For any subset $I$ of
$\left\{  1,2,\ldots,n\right\}  $, we let $\widetilde{I}$ denote the
complement $\left\{  1,2,\ldots,n\right\}  \setminus I$ of $I$.

Let $A=\left(  a_{i,j}\right)  _{1\leq i\leq n,\ 1\leq j\leq n}$ and
$B=\left(  b_{i,j}\right)  _{1\leq i\leq n,\ 1\leq j\leq n}$ be two $n\times
n$-matrices. Let $P$ and $Q$ be two subsets of $\left\{  1,2,\ldots,n\right\}
$ such that $\left\vert P\right\vert =\left\vert Q\right\vert $. Then,%
\begin{align*}
&  \sum_{\substack{\sigma\in S_{n};\\\sigma\left(  P\right)  =Q}}\left(
-1\right)  ^{\sigma}\left(  \prod_{i\in P}a_{i,\sigma\left(  i\right)
}\right)  \left(  \prod_{i\in\widetilde{P}}b_{i,\sigma\left(  i\right)
}\right) \\
&  =\left(  -1\right)  ^{\sum P+\sum Q}\det\left(  \operatorname*{sub}%
\nolimits_{w\left(  P\right)  }^{w\left(  Q\right)  }A\right)  \det\left(
\operatorname*{sub}\nolimits_{w\left(  \widetilde{P}\right)  }^{w\left(
\widetilde{Q}\right)  }B\right)  .
\end{align*}

\end{lemma}

The proof of Lemma \ref{lem.det.laplace-multi.Apq} is similar (in its spirit)
to the proof of Lemma \ref{lem.laplace.Apq}, but it requires a lot more
bookkeeping (if one wants to make it rigorous). This proof shall be given in
the solution to Exercise \ref{exe.det.laplace-multi}:

\begin{exercise}
\label{exe.det.laplace-multi}Prove Lemma \ref{lem.det.laplace-multi.Apq} and
Theorem \ref{thm.det.laplace-multi}.

[\textbf{Hint:} First, prove Lemma \ref{lem.det.laplace-multi.Apq} in the case
when $P=\left\{  1,2,\ldots,k\right\}  $ for some $k\in\left\{  0,1,\ldots
,n\right\}  $; in order to do so, use the bijection from Exercise
\ref{exe.Ialbe} \textbf{(c)} (applied to $I=Q$). Then, derive the general case
of Lemma \ref{lem.det.laplace-multi.Apq} by permuting the rows of the
matrices. Finally, prove Theorem \ref{thm.det.laplace-multi}.]
\end{exercise}

The following exercise generalizes Proposition \ref{prop.laplace.0} in the
same way as Theorem \ref{thm.det.laplace-multi} generalizes Theorem
\ref{thm.laplace.gen}:

\begin{exercise}
\label{exe.det.laplace-multi.0}Let $n\in\mathbb{N}$. Let $A\in\mathbb{K}%
^{n\times n}$. For any subset $I$ of $\left\{  1,2,\ldots,n\right\}  $, we let
$\widetilde{I}$ denote the complement $\left\{  1,2,\ldots,n\right\}
\setminus I$ of $I$. Let $R$ be a subset of $\left\{  1,2,\ldots,n\right\}  $.
Prove the following:

\textbf{(a)} For every subset $P$ of $\left\{  1,2,\ldots,n\right\}  $
satisfying $\left\vert P\right\vert =\left\vert R\right\vert $ and $P\neq R$,
we have%
\[
0=\sum_{\substack{Q\subseteq\left\{  1,2,\ldots,n\right\}  ;\\\left\vert
Q\right\vert =\left\vert P\right\vert }}\left(  -1\right)  ^{\sum P+\sum
Q}\det\left(  \operatorname*{sub}\nolimits_{w\left(  R\right)  }^{w\left(
Q\right)  }A\right)  \det\left(  \operatorname*{sub}\nolimits_{w\left(
\widetilde{P}\right)  }^{w\left(  \widetilde{Q}\right)  }A\right)  .
\]

\textbf{(b)} For every subset $Q$ of $\left\{  1,2,\ldots,n\right\}  $
satisfying $\left\vert Q\right\vert =\left\vert R\right\vert $ and $Q\neq R$,
we have%
\[
0=\sum_{\substack{P\subseteq\left\{  1,2,\ldots,n\right\}  ;\\\left\vert
P\right\vert =\left\vert Q\right\vert }}\left(  -1\right)  ^{\sum P+\sum
Q}\det\left(  \operatorname*{sub}\nolimits_{w\left(  P\right)  }^{w\left(
R\right)  }A\right)  \det\left(  \operatorname*{sub}\nolimits_{w\left(
\widetilde{P}\right)  }^{w\left(  \widetilde{Q}\right)  }A\right)  .
\]

\end{exercise}

Exercise \ref{exe.det.laplace-multi.0} can be generalized to non-square matrices:

\begin{exercise}
\label{exe.det.laplace-multi.0r}Let $n\in\mathbb{N}$ and $m\in\mathbb{N}$. For
any subset $I$ of $\left\{  1,2,\ldots,n\right\}  $, we let $\widetilde{I}$
denote the complement $\left\{  1,2,\ldots,n\right\}  \setminus I$ of $I$. Let
$J$ and $K$ be two subsets of $\left\{  1,2,\ldots,m\right\}  $ satisfying
$\left\vert J\right\vert +\left\vert K\right\vert =n$ and $J\cap
K\neq\varnothing$. Prove the following:

\textbf{(a)} For every $A\in\mathbb{K}^{m\times n}$, we have%
\[
0=\sum_{\substack{Q\subseteq\left\{  1,2,\ldots,n\right\}  ;\\\left\vert
Q\right\vert =\left\vert J\right\vert }}\left(  -1\right)  ^{\sum Q}%
\det\left(  \operatorname*{sub}\nolimits_{w\left(  J\right)  }^{w\left(
Q\right)  }A\right)  \det\left(  \operatorname*{sub}\nolimits_{w\left(
K\right)  }^{w\left(  \widetilde{Q}\right)  }A\right)  .
\]

\textbf{(b)} For every $A\in\mathbb{K}^{n\times m}$, we have%
\[
0=\sum_{\substack{P\subseteq\left\{  1,2,\ldots,n\right\}  ;\\\left\vert
P\right\vert =\left\vert J\right\vert }}\left(  -1\right)  ^{\sum P}%
\det\left(  \operatorname*{sub}\nolimits_{w\left(  P\right)  }^{w\left(
J\right)  }A\right)  \det\left(  \operatorname*{sub}\nolimits_{w\left(
\widetilde{P}\right)  }^{w\left(  K\right)  }A\right)  .
\]

\end{exercise}

(Exercise \ref{exe.det.laplace-multi.0} can be obtained as a particular case
of Exercise \ref{exe.det.laplace-multi.0r}; we leave the details to the reader.)

The following exercise gives a first application of Theorem
\ref{thm.det.laplace-multi} (though it can also be solved with more elementary methods):

\begin{exercise}
\label{exe.det.blocktria-twisted}Let $n\in\mathbb{N}$. Let $P$ and $Q$ be two
subsets of $\left\{  1,2,\ldots,n\right\}  $. Let $A=\left(  a_{i,j}\right)
_{1\leq i\leq n,\ 1\leq j\leq n}\in\mathbb{K}^{n\times n}$ be an $n\times
n$-matrix such that%
\begin{equation}
\text{every }i\in P\text{ and }j\in Q\text{ satisfy }a_{i,j}=0.
\label{eq.exe.det.blocktria-twisted.ass}%
\end{equation}
For any subset $I$ of $\left\{  1,2,\ldots,n\right\}  $, we let $\widetilde{I}%
$ denote the complement $\left\{  1,2,\ldots,n\right\}  \setminus I$ of $I$.
(For instance, if $n=4$ and $I=\left\{  1,4\right\}  $, then $\widetilde{I}%
=\left\{  2,3\right\}  $.)

Prove the following:

\textbf{(a)} If $\left\vert P\right\vert +\left\vert Q\right\vert >n$, then
$\det A=0$.

\textbf{(b)} If $\left\vert P\right\vert +\left\vert Q\right\vert =n$, then
\[
\det A=\left(  -1\right)  ^{\sum P+\sum\widetilde{Q}}\det\left(
\operatorname*{sub}\nolimits_{w\left(  P\right)  }^{w\left(  \widetilde{Q}%
\right)  }A\right)  \det\left(  \operatorname*{sub}\nolimits_{w\left(
\widetilde{P}\right)  }^{w\left(  Q\right)  }A\right)  .
\]

\end{exercise}

\begin{example}
\label{exa.det.blocktria-twisted.exas}\textbf{(a)} Applying Exercise
\ref{exe.det.blocktria-twisted} \textbf{(a)} to $n=5$, $P=\left\{
1,3,5\right\}  $ and $Q=\left\{  2,3,4\right\}  $, we see that%
\[
\det\left(
\begin{array}
[c]{ccccc}%
a & 0 & 0 & 0 & b\\
c & d & e & f & g\\
h & 0 & 0 & 0 & i\\
j & k & \ell & m & n\\
o & 0 & 0 & 0 & p
\end{array}
\right)  =0
\]
(don't mistake the letter \textquotedblleft$o$\textquotedblright\ for a zero)
for any $a,b,c,d,e,f,g,h,i,j,k,\ell,m,n,o,p\in\mathbb{K}$ (since the $n\times
n$-matrices $A=\left(  a_{i,j}\right)  _{1\leq i\leq n,\ 1\leq j\leq n}$
satisfying (\ref{eq.exe.det.blocktria-twisted.ass}) are precisely the matrices
of the form $\left(
\begin{array}
[c]{ccccc}%
a & 0 & 0 & 0 & b\\
c & d & e & f & g\\
h & 0 & 0 & 0 & i\\
j & k & \ell & m & n\\
o & 0 & 0 & 0 & p
\end{array}
\right)  $).

\textbf{(b)} Applying Exercise \ref{exe.det.blocktria-twisted} \textbf{(a)} to
$n=5$, $P=\left\{  2,3,4\right\}  $ and $Q=\left\{  2,3,4\right\}  $, we see
that%
\[
\det\left(
\begin{array}
[c]{ccccc}%
a & b & c & d & e\\
f & 0 & 0 & 0 & g\\
h & 0 & 0 & 0 & i\\
j & 0 & 0 & 0 & k\\
\ell & m & n & o & p
\end{array}
\right)  =0
\]
(don't mistake the letter \textquotedblleft$o$\textquotedblright\ for a zero)
for any $a,b,c,d,e,f,g,h,i,j,k,\ell,m,n,o,p\in\mathbb{K}$. This is precisely
the claim of Exercise \ref{exe.ps4.5} \textbf{(b)}.

\textbf{(c)} Applying Exercise \ref{exe.det.blocktria-twisted} \textbf{(b)} to
$n=4$, $P=\left\{  2,3\right\}  $ and $Q=\left\{  2,3\right\}  $, we see that%
\[
\det\left(
\begin{array}
[c]{cccc}%
a & b & c & d\\
\ell & 0 & 0 & e\\
k & 0 & 0 & f\\
j & i & h & g
\end{array}
\right)  =\det\left(
\begin{array}
[c]{cc}%
\ell & e\\
k & f
\end{array}
\right)  \cdot\det\left(
\begin{array}
[c]{cc}%
b & c\\
i & h
\end{array}
\right)
\]
for all $a,b,c,d,e,f,g,h,i,j,k,\ell\in\mathbb{K}$. This solves Exercise
\ref{exe.ps4.5} \textbf{(a)}.

\textbf{(d)} It is not hard to derive Exercise \ref{exe.block2x2.tridet} by
applying Exercise \ref{exe.det.blocktria-twisted} \textbf{(b)} to $n+m$,
$\left(
\begin{array}
[c]{cc}%
A & 0_{n\times m}\\
C & D
\end{array}
\right)  $, $\left\{  1,2,\ldots,n\right\}  $ and $\left\{  n+1,n+2,\ldots
,n+m\right\}  $ instead of $n$, $A$, $P$ and $Q$. Similarly, Exercise
\ref{exe.block2x2.tridet.transposed} can be derived from Exercise
\ref{exe.det.blocktria-twisted} \textbf{(b)} as well.
\end{example}

\subsection{$\det\left(  A+B\right)  $}

As Theorem \ref{thm.det(AB)} shows, the determinant of the product $AB$ of two
square matrices can be easily and neatly expressed through the determinants of
$A$ and $B$. In contrast, the determinant of a sum $A+B$ of two square
matrices cannot be expressed in such a way\footnote{It is easy to find two
$2\times2$-matrices $A_{1}$ and $B_{1}$ and two other $2\times2$-matrices
$A_{2}$ and $B_{2}$ such that $\det\left(  A_{1}\right)  =\det\left(
A_{2}\right)  $ and $\det\left(  B_{1}\right)  =\det\left(  B_{2}\right)  $
but $\det\left(  A_{1}+B_{1}\right)  \neq\det\left(  A_{2}+B_{2}\right)  $.
This shows that $\det\left(  A+B\right)  $ cannot generally be computed from
$\det A$ and $\det B$.}. There is, however, a formula for $\det\left(
A+B\right)  $ in terms of the determinants of submatrices of $A$ and $B$.
While it is rather unwieldy (a far cry from the elegance of Theorem
\ref{thm.det(AB)}), it is nevertheless useful sometimes; let us now show it:

\begin{theorem}
\label{thm.det(A+B)}Let $n\in\mathbb{N}$. For any subset $I$ of $\left\{
1,2,\ldots,n\right\}  $, we let $\widetilde{I}$ denote the complement
$\left\{  1,2,\ldots,n\right\}  \setminus I$ of $I$. (For instance, if $n=4$
and $I=\left\{  1,4\right\}  $, then $\widetilde{I}=\left\{  2,3\right\}  $.)
Let us use the notations introduced in Definition \ref{def.submatrix} and in
Definition \ref{def.sect.laplace.notations}.

Let $A$ and $B$ be two $n\times n$-matrices. Then,%
\[
\det\left(  A+B\right)  =\sum_{P\subseteq\left\{  1,2,\ldots,n\right\}  }%
\sum_{\substack{Q\subseteq\left\{  1,2,\ldots,n\right\}  ;\\\left\vert
P\right\vert =\left\vert Q\right\vert }}\left(  -1\right)  ^{\sum P+\sum
Q}\det\left(  \operatorname*{sub}\nolimits_{w\left(  P\right)  }^{w\left(
Q\right)  }A\right)  \det\left(  \operatorname*{sub}\nolimits_{w\left(
\widetilde{P}\right)  }^{w\left(  \widetilde{Q}\right)  }B\right)  .
\]

\end{theorem}

\begin{example}
For this example, set $n=2$, $A=\left(  a_{i,j}\right)  _{1\leq i\leq2,\ 1\leq
j\leq2}$ and $B=\left(  b_{i,j}\right)  _{1\leq i\leq2,\ 1\leq j\leq2}$. Then,
Theorem \ref{thm.det(A+B)} says that%
\begin{align*}
&  \det\left(  A+B\right) \\
&  =\sum_{P\subseteq\left\{  1,2\right\}  }\sum_{\substack{Q\subseteq\left\{
1,2\right\}  ;\\\left\vert P\right\vert =\left\vert Q\right\vert }}\left(
-1\right)  ^{\sum P+\sum Q}\det\left(  \operatorname*{sub}\nolimits_{w\left(
P\right)  }^{w\left(  Q\right)  }A\right)  \det\left(  \operatorname*{sub}%
\nolimits_{w\left(  \widetilde{P}\right)  }^{w\left(  \widetilde{Q}\right)
}B\right) \\
&  =\left(  -1\right)  ^{\sum\varnothing+\sum\varnothing}\det
\underbrace{\left(  \operatorname*{sub}\nolimits_{{}}^{{}}A\right)
}_{\substack{\text{this is the}\\0\times0\text{-matrix}}}\det\left(
\operatorname*{sub}\nolimits_{1,2}^{1,2}B\right) \\
&  \ \ \ \ \ \ \ \ \ \ +\left(  -1\right)  ^{\sum\left\{  1\right\}
+\sum\left\{  1\right\}  }\det\left(  \operatorname*{sub}\nolimits_{1}%
^{1}A\right)  \det\left(  \operatorname*{sub}\nolimits_{2}^{2}B\right) \\
&  \ \ \ \ \ \ \ \ \ \ +\left(  -1\right)  ^{\sum\left\{  1\right\}
+\sum\left\{  2\right\}  }\det\left(  \operatorname*{sub}\nolimits_{1}%
^{2}A\right)  \det\left(  \operatorname*{sub}\nolimits_{2}^{1}B\right) \\
&  \ \ \ \ \ \ \ \ \ \ +\left(  -1\right)  ^{\sum\left\{  2\right\}
+\sum\left\{  1\right\}  }\det\left(  \operatorname*{sub}\nolimits_{2}%
^{1}A\right)  \det\left(  \operatorname*{sub}\nolimits_{1}^{2}B\right) \\
&  \ \ \ \ \ \ \ \ \ \ +\left(  -1\right)  ^{\sum\left\{  2\right\}
+\sum\left\{  2\right\}  }\det\left(  \operatorname*{sub}\nolimits_{2}%
^{2}A\right)  \det\left(  \operatorname*{sub}\nolimits_{1}^{1}B\right) \\
&  \ \ \ \ \ \ \ \ \ \ +\left(  -1\right)  ^{\sum\left\{  1,2\right\}
+\sum\left\{  1,2\right\}  }\det\left(  \operatorname*{sub}\nolimits_{1,2}%
^{1,2}A\right)  \det\underbrace{\left(  \operatorname*{sub}\nolimits_{{}}^{{}%
}B\right)  }_{\substack{\text{this is the}\\0\times0\text{-matrix}}}\\
&  =\underbrace{\det\left(  \text{the }0\times0\text{-matrix}\right)  }%
_{=1}\det\left(
\begin{array}
[c]{cc}%
b_{1,1} & b_{1,2}\\
b_{2,1} & b_{2,2}%
\end{array}
\right)  +\det\left(
\begin{array}
[c]{c}%
a_{1,1}%
\end{array}
\right)  \det\left(
\begin{array}
[c]{c}%
b_{2,2}%
\end{array}
\right) \\
&  \ \ \ \ \ \ \ \ \ \ -\det\left(
\begin{array}
[c]{c}%
a_{1,2}%
\end{array}
\right)  \det\left(
\begin{array}
[c]{c}%
b_{2,1}%
\end{array}
\right)  -\det\left(
\begin{array}
[c]{c}%
a_{2,1}%
\end{array}
\right)  \det\left(
\begin{array}
[c]{c}%
b_{1,2}%
\end{array}
\right) \\
&  \ \ \ \ \ \ \ \ \ \ +\det\left(
\begin{array}
[c]{c}%
a_{2,2}%
\end{array}
\right)  \det\left(
\begin{array}
[c]{c}%
b_{1,1}%
\end{array}
\right)  +\det\left(
\begin{array}
[c]{cc}%
a_{1,1} & a_{1,2}\\
a_{2,1} & a_{2,2}%
\end{array}
\right)  \underbrace{\det\left(  \text{the }0\times0\text{-matrix}\right)
}_{=1}\\
&  =\det\left(
\begin{array}
[c]{cc}%
b_{1,1} & b_{1,2}\\
b_{2,1} & b_{2,2}%
\end{array}
\right)  +a_{1,1}b_{2,2}-a_{1,2}b_{2,1}-a_{2,1}b_{1,2}+a_{2,2}b_{1,1}%
+\det\left(
\begin{array}
[c]{cc}%
a_{1,1} & a_{1,2}\\
a_{2,1} & a_{2,2}%
\end{array}
\right)  .
\end{align*}

\end{example}

\begin{exercise}
\label{exe.det(A+B)}Prove Theorem \ref{thm.det(A+B)}.

[\textbf{Hint:} Use Lemma \ref{lem.det.laplace-multi.Apq}.]
\end{exercise}

Theorem \ref{thm.det(A+B)} takes a simpler form in the particular case when
the matrix $B$ is diagonal (i.e., has all entries outside of its diagonal
equal to $0$):

\begin{corollary}
\label{cor.det(A+D)}Let $n\in\mathbb{N}$. For every two objects $i$ and $j$,
define $\delta_{i,j}\in\mathbb{K}$ by $\delta_{i,j}=%
\begin{cases}
1, & \text{if }i=j;\\
0, & \text{if }i\neq j
\end{cases}
$.

Let $A$ be an $n\times n$-matrix. Let $d_{1},d_{2},\ldots,d_{n}$ be $n$
elements of $\mathbb{K}$. Let $D$ be the $n\times n$-matrix $\left(
d_{i}\delta_{i,j}\right)  _{1\leq i\leq n,\ 1\leq j\leq n}$. Then,%
\[
\det\left(  A+D\right)  =\sum_{P\subseteq\left\{  1,2,\ldots,n\right\}  }%
\det\left(  \operatorname*{sub}\nolimits_{w\left(  P\right)  }^{w\left(
P\right)  }A\right)  \prod_{i\in\left\{  1,2,\ldots,n\right\}  \setminus
P}d_{i}.
\]

\end{corollary}

This corollary can easily be derived from Theorem \ref{thm.det(A+B)} using the
following fact:

\begin{lemma}
\label{lem.diag.minors}Let $n\in\mathbb{N}$. For every two objects $i$ and
$j$, define $\delta_{i,j}\in\mathbb{K}$ by $\delta_{i,j}=%
\begin{cases}
1, & \text{if }i=j;\\
0, & \text{if }i\neq j
\end{cases}
$.

Let $d_{1},d_{2},\ldots,d_{n}$ be $n$ elements of $\mathbb{K}$. Let $D$ be the
$n\times n$-matrix $\left(  d_{i}\delta_{i,j}\right)  _{1\leq i\leq n,\ 1\leq
j\leq n}$. Let $P$ and $Q$ be two subsets of $\left\{  1,2,\ldots,n\right\}  $
such that $\left\vert P\right\vert =\left\vert Q\right\vert $. Then,%
\[
\det\left(  \operatorname*{sub}\nolimits_{w\left(  P\right)  }^{w\left(
Q\right)  }D\right)  =\delta_{P,Q}\prod_{i\in P}d_{i}.
\]

\end{lemma}

Proving Corollary \ref{cor.det(A+D)} and Lemma \ref{lem.diag.minors} in detail
is part of Exercise \ref{exe.det(A+B).diag} further below.

A particular case of Corollary \ref{cor.det(A+D)} is the following fact:

\begin{corollary}
\label{cor.det(A+X)}Let $n\in\mathbb{N}$. Let $A$ be an $n\times n$-matrix.
Let $x\in\mathbb{K}$. Then,%
\begin{align}
\det\left(  A+xI_{n}\right)   &  =\sum_{P\subseteq\left\{  1,2,\ldots
,n\right\}  }\det\left(  \operatorname*{sub}\nolimits_{w\left(  P\right)
}^{w\left(  P\right)  }A\right)  x^{n-\left\vert P\right\vert }%
\label{eq.cor.det(A+X).1}\\
&  =\sum_{k=0}^{n}\left(  \sum_{\substack{P\subseteq\left\{  1,2,\ldots
,n\right\}  ;\\\left\vert P\right\vert =n-k}}\det\left(  \operatorname*{sub}%
\nolimits_{w\left(  P\right)  }^{w\left(  P\right)  }A\right)  \right)  x^{k}.
\label{eq.cor.det(A+X).2}%
\end{align}

\end{corollary}

\begin{exercise}
\label{exe.det(A+B).diag}Prove Corollary \ref{cor.det(A+D)}, Lemma
\ref{lem.diag.minors} and Corollary \ref{cor.det(A+X)}.
\end{exercise}

\begin{remark}
\label{rmk.charpol}Let $n\in\mathbb{N}$. Let $A$ be an $n\times n$-matrix over
the commutative ring $\mathbb{K}$. Consider the commutative ring
$\mathbb{K}\left[  X\right]  $ of polynomials in the indeterminate $X$ over
$\mathbb{K}$ (that is, polynomials in the indeterminate $X$ with coefficients
lying in $\mathbb{K}$). We can then regard $A$ as a matrix over the ring
$\mathbb{K}\left[  X\right]  $ as well (because every element of $\mathbb{K}$
can be viewed as a constant polynomial in $\mathbb{K}\left[  X\right]  $).

Consider the $n\times n$-matrix $A+XI_{n}$ over the commutative ring
$\mathbb{K}\left[  X\right]  $. (For example, if $n=2$ and $A=\left(
\begin{array}
[c]{cc}%
a & b\\
c & d
\end{array}
\right)  $, then $A+XI_{n}=\left(
\begin{array}
[c]{cc}%
a & b\\
c & d
\end{array}
\right)  +X\left(
\begin{array}
[c]{cc}%
1 & 0\\
0 & 1
\end{array}
\right)  =\left(
\begin{array}
[c]{cc}%
a+X & b\\
c & d+X
\end{array}
\right)  $. In general, the matrix $A+XI_{n}$ is obtained from $A$ by adding
an $X$ to each diagonal entry.)

The determinant $\det\left(  A+XI_{n}\right)  $ is a polynomial in
$\mathbb{K}\left[  X\right]  $. (For instance, for $n=2$ and $A=\left(
\begin{array}
[c]{cc}%
a & b\\
c & d
\end{array}
\right)  $, we have%
\begin{align*}
\det\left(  A+XI_{n}\right)   &  =\det\left(
\begin{array}
[c]{cc}%
a+X & b\\
c & d+X
\end{array}
\right)  =\left(  a+X\right)  \left(  d+X\right)  -bc\\
&  =X^{2}+\left(  a+d\right)  X+\left(  ad-bc\right)  .
\end{align*}
)

This polynomial $\det\left(  A+XI_{n}\right)  $ is a highly important object;
it is a close relative of what is called the \textit{characteristic
polynomial} of $A$. (More precisely, the characteristic polynomial of $A$ is
either $\det\left(  XI_{n}-A\right)  $ or $\det\left(  A-XI_{n}\right)  $,
depending on the conventions that one is using; thus, the polynomial
$\det\left(  A+XI_{n}\right)  $ is either the characteristic polynomial of
$-A$ or $\left(  -1\right)  ^{n}$ times this characteristic polynomial.) For
more about the characteristic polynomial, see \cite[Section 4.5]{Artin} or
\cite[Chapter Five, Section II, \S 3]{Hefferon} (or various other texts on
linear algebra).

Using Corollary \ref{cor.det(A+X)}, we can explicitly compute the coefficients
of the polynomial $\det\left(  A+XI_{n}\right)  $. In fact,
(\ref{eq.cor.det(A+X).2}) (applied to $\mathbb{K}\left[  X\right]  $ and $X$
instead of $\mathbb{K}$ and $x$) yields%
\[
\det\left(  A+XI_{n}\right)  =\sum_{k=0}^{n}\left(  \sum_{\substack{P\subseteq
\left\{  1,2,\ldots,n\right\}  ;\\\left\vert P\right\vert =n-k}}\det\left(
\operatorname*{sub}\nolimits_{w\left(  P\right)  }^{w\left(  P\right)
}A\right)  \right)  X^{k}.
\]
Hence, for every $k\in\left\{  0,1,\ldots,n\right\}  $, the coefficient of
$X^{k}$ in the polynomial $\det\left(  A+XI_{n}\right)  $ is
\[
\sum_{\substack{P\subseteq\left\{  1,2,\ldots,n\right\}  ;\\\left\vert
P\right\vert =n-k}}\det\left(  \operatorname*{sub}\nolimits_{w\left(
P\right)  }^{w\left(  P\right)  }A\right)  .
\]
In particular:

\begin{itemize}
\item The coefficient of $X^{n}$ in the polynomial $\det\left(  A+XI_{n}%
\right)  $ is%
\begin{align*}
&  \sum_{\substack{P\subseteq\left\{  1,2,\ldots,n\right\}  ;\\\left\vert
P\right\vert =n-n}}\det\left(  \operatorname*{sub}\nolimits_{w\left(
P\right)  }^{w\left(  P\right)  }A\right) \\
&  =\det\underbrace{\left(  \operatorname*{sub}\nolimits_{w\left(
\varnothing\right)  }^{w\left(  \varnothing\right)  }A\right)  }_{=\left(
\text{the }0\times0\text{-matrix}\right)  }\\
&  \ \ \ \ \ \ \ \ \ \ \left(
\begin{array}
[c]{c}%
\text{since the only subset }P\text{ of }\left\{  1,2,\ldots,n\right\} \\
\text{satisfying }\left\vert P\right\vert =n-n\text{ is the empty set
}\varnothing
\end{array}
\right) \\
&  =\det\left(  \text{the }0\times0\text{-matrix}\right)  =1.
\end{align*}

\item The coefficient of $X^{0}$ in the polynomial $\det\left(  A+XI_{n}%
\right)  $ is%
\begin{align*}
&  \sum_{\substack{P\subseteq\left\{  1,2,\ldots,n\right\}  ;\\\left\vert
P\right\vert =n-0}}\det\left(  \operatorname*{sub}\nolimits_{w\left(
P\right)  }^{w\left(  P\right)  }A\right) \\
&  =\det\underbrace{\left(  \operatorname*{sub}\nolimits_{w\left(  \left\{
1,2,\ldots,n\right\}  \right)  }^{w\left(  \left\{  1,2,\ldots,n\right\}
\right)  }A\right)  }_{=\operatorname*{sub}\nolimits_{1,2,\ldots
,n}^{1,2,\ldots,n}A=A}\\
&  \ \ \ \ \ \ \ \ \ \ \left(
\begin{array}
[c]{c}%
\text{since the only subset }P\text{ of }\left\{  1,2,\ldots,n\right\} \\
\text{satisfying }\left\vert P\right\vert =n-0\text{ is the set }\left\{
1,2,\ldots,n\right\}
\end{array}
\right) \\
&  =\det A.
\end{align*}

\item Write the matrix $A$ as $A=\left(  a_{i,j}\right)  _{1\leq i\leq
n,\ 1\leq j\leq n}$. Then, the coefficient of $X^{n-1}$ in the polynomial
$\det\left(  A+XI_{n}\right)  $ is%
\begin{align*}
&  \sum_{\substack{P\subseteq\left\{  1,2,\ldots,n\right\}  ;\\\left\vert
P\right\vert =n-\left(  n-1\right)  }}\det\left(  \operatorname*{sub}%
\nolimits_{w\left(  P\right)  }^{w\left(  P\right)  }A\right) \\
&  =\sum_{\substack{P\subseteq\left\{  1,2,\ldots,n\right\}  ;\\\left\vert
P\right\vert =1}}\det\left(  \operatorname*{sub}\nolimits_{w\left(  P\right)
}^{w\left(  P\right)  }A\right)  =\sum_{k=1}^{n}\det\underbrace{\left(
\operatorname*{sub}\nolimits_{w\left(  \left\{  k\right\}  \right)
}^{w\left(  \left\{  k\right\}  \right)  }A\right)  }%
_{\substack{=\operatorname*{sub}\nolimits_{k}^{k}A=\left(
\begin{array}
[c]{c}%
a_{k,k}%
\end{array}
\right)  \\\text{(this is a }1\times1\text{-matrix)}}}\\
&  \ \ \ \ \ \ \ \ \ \ \left(
\begin{array}
[c]{c}%
\text{since the subsets }P\text{ of }\left\{  1,2,\ldots,n\right\} \\
\text{satisfying }\left\vert P\right\vert =1\text{ are the sets }\left\{
1\right\}  ,\left\{  2\right\}  ,\ldots,\left\{  n\right\}
\end{array}
\right) \\
&  =\sum_{k=1}^{n}\underbrace{\det\left(
\begin{array}
[c]{c}%
a_{k,k}%
\end{array}
\right)  }_{=a_{k,k}}=\sum_{k=1}^{n}a_{k,k}.
\end{align*}
In other words, this coefficient is the sum of all diagonal entries of $A$.
This sum is called the \textit{trace} of $A$, and is denoted by
$\operatorname*{Tr}A$.
\end{itemize}
\end{remark}

\subsection{Some alternating-sum formulas}

The next few exercises don't all involve determinants; what they have in
common is that they contain alternating sums (i.e., sums where the addend
contains a power of $-1$). \Needspace{8cm}

\begin{exercise}
\label{exe.noncomm.polarization}For every $n\in\mathbb{N}$, let $\left[
n\right]  $ denote the set $\left\{  1,2,\ldots,n\right\}  $.

Let $\mathbb{L}$ be a noncommutative ring. (Keep in mind that our definition
of a \textquotedblleft noncommutative ring\textquotedblright\ includes all
commutative rings.)

Let $n\in\mathbb{N}$. The summation sign $\sum_{I\subseteq\left[  n\right]  }$
shall mean $\sum_{I\in\mathcal{P}\left(  \left[  n\right]  \right)  }$, where
$\mathcal{P}\left(  \left[  n\right]  \right)  $ denotes the powerset of
$\left[  n\right]  $.

Let $v_{1},v_{2},\ldots,v_{n}$ be $n$ elements of $\mathbb{L}$.

\textbf{(a)} Prove that%
\[
\sum_{I\subseteq\left[  n\right]  }\left(  -1\right)  ^{n-\left\vert
I\right\vert }\left(  \sum_{i\in I}v_{i}\right)  ^{m}=\sum
_{\substack{f:\left[  m\right]  \rightarrow\left[  n\right]  ;\\f\text{ is
surjective}}}v_{f\left(  1\right)  }v_{f\left(  2\right)  }\cdots v_{f\left(
m\right)  }%
\]
for each $m\in\mathbb{N}$.

\textbf{(b)} Prove that%
\[
\sum_{I\subseteq\left[  n\right]  }\left(  -1\right)  ^{n-\left\vert
I\right\vert }\left(  \sum_{i\in I}v_{i}\right)  ^{m}=0
\]
for each $m\in\left\{  0,1,\ldots,n-1\right\}  $.

\textbf{(c)} Prove that%
\[
\sum_{I\subseteq\left[  n\right]  }\left(  -1\right)  ^{n-\left\vert
I\right\vert }\left(  \sum_{i\in I}v_{i}\right)  ^{n}=\sum_{\sigma\in S_{n}%
}v_{\sigma\left(  1\right)  }v_{\sigma\left(  2\right)  }\cdots v_{\sigma
\left(  n\right)  }.
\]

\textbf{(d)} Now, assume that $\mathbb{L}$ is a \textbf{commutative} ring.
Prove that%
\[
\sum_{I\subseteq\left[  n\right]  }\left(  -1\right)  ^{n-\left\vert
I\right\vert }\left(  \sum_{i\in I}v_{i}\right)  ^{n}=n!\cdot v_{1}v_{2}\cdots
v_{n}.
\]

[\textbf{Hint:} First, generalize Lemma \ref{lem.prodrule2} to the case of a
noncommutative ring $\mathbb{K}$.]
\end{exercise}

\begin{exercise}
\label{exe.noncomm.polarization2}For every $n\in\mathbb{N}$, let $\left[
n\right]  $ denote the set $\left\{  1,2,\ldots,n\right\}  $.

Let $\mathbb{L}$ be a noncommutative ring. (Keep in mind that our definition
of a \textquotedblleft noncommutative ring\textquotedblright\ includes all
commutative rings.)

Let $n\in\mathbb{N}$. The summation sign $\sum_{I\subseteq\left[  n\right]  }$
shall mean $\sum_{I\in\mathcal{P}\left(  \left[  n\right]  \right)  }$, where
$\mathcal{P}\left(  \left[  n\right]  \right)  $ denotes the powerset of
$\left[  n\right]  $.

Let $v_{1},v_{2},\ldots,v_{n}$ be $n$ elements of $\mathbb{L}$.

\textbf{(a)} Prove that%
\[
\sum_{I\subseteq\left[  n\right]  }\left(  -1\right)  ^{n-\left\vert
I\right\vert }\left(  w+\sum_{i\in I}v_{i}\right)  ^{m}=0
\]
for each $m\in\left\{  0,1,\ldots,n-1\right\}  $ and $w\in\mathbb{L}$.

\textbf{(b)} Prove that
\[
\sum_{I\subseteq\left[  n\right]  }\left(  -1\right)  ^{n-\left\vert
I\right\vert }\left(  w+\sum_{i\in I}v_{i}\right)  ^{n}=\sum_{\sigma\in S_{n}%
}v_{\sigma\left(  1\right)  }v_{\sigma\left(  2\right)  }\cdots v_{\sigma
\left(  n\right)  }%
\]
for every $w\in\mathbb{L}$.

\textbf{(c)} Prove that%
\[
\sum_{I\subseteq\left[  n\right]  }\left(  -1\right)  ^{n-\left\vert
I\right\vert }\left(  \sum_{i\in I}v_{i}-\sum_{i\in\left[  n\right]  \setminus
I}v_{i}\right)  ^{m}=0
\]
for each $m\in\left\{  0,1,\ldots,n-1\right\}  $.

\textbf{(d)} Prove that
\[
\sum_{I\subseteq\left[  n\right]  }\left(  -1\right)  ^{n-\left\vert
I\right\vert }\left(  \sum_{i\in I}v_{i}-\sum_{i\in\left[  n\right]  \setminus
I}v_{i}\right)  ^{n}=2^{n}\sum_{\sigma\in S_{n}}v_{\sigma\left(  1\right)
}v_{\sigma\left(  2\right)  }\cdots v_{\sigma\left(  n\right)  }.
\]

\end{exercise}

Note that parts \textbf{(b)} and \textbf{(c)} of Exercise
\ref{exe.noncomm.polarization} are special cases of parts \textbf{(a)} and
\textbf{(b)} of Exercise \ref{exe.noncomm.polarization2} (obtained by setting
$w=0$).

The following exercise generalizes Exercise \ref{exe.noncomm.polarization}:

\begin{exercise}
\label{exe.det.sumdets1}For every $n\in\mathbb{N}$, let $\left[  n\right]  $
denote the set $\left\{  1,2,\ldots,n\right\}  $.

Let $\mathbb{L}$ be a noncommutative ring. (Keep in mind that our definition
of a \textquotedblleft noncommutative ring\textquotedblright\ includes all
commutative rings.)

Let $G$ be a finite set. Let $H$ be a subset of $G$. Let $n\in\mathbb{N}$.

For each $i\in G$ and $j\in\left[  n\right]  $, let $b_{i,j}$ be an element of
$\mathbb{L}$. For each $j\in\left[  n\right]  $ and each subset $I$ of $G$, we
define an element $b_{I,j}\in\mathbb{L}$ by $b_{I,j}=\sum_{i\in I}b_{i,j}$.

\textbf{(a)} Prove that%
\[
\sum_{\substack{I\subseteq G;\\H\subseteq I}}\left(  -1\right)  ^{\left\vert
I\right\vert }b_{I,1}b_{I,2}\cdots b_{I,n}=\left(  -1\right)  ^{\left\vert
G\right\vert }\sum_{\substack{f:\left[  n\right]  \rightarrow G;\\G\setminus
H\subseteq f\left(  \left[  n\right]  \right)  }}b_{f\left(  1\right)
,1}b_{f\left(  2\right)  ,2}\cdots b_{f\left(  n\right)  ,n}.
\]

\textbf{(b)} If $n<\left\vert G\setminus H\right\vert $, then prove that%
\[
\sum_{\substack{I\subseteq G;\\H\subseteq I}}\left(  -1\right)  ^{\left\vert
I\right\vert }b_{I,1}b_{I,2}\cdots b_{I,n}=0.
\]

\textbf{(c)} If $n=\left\vert G\right\vert $, then prove that%
\[
\sum_{I\subseteq G}\left(  -1\right)  ^{\left\vert G\setminus I\right\vert
}b_{I,1}b_{I,2}\cdots b_{I,n}=\sum_{\substack{f:\left[  n\right]  \rightarrow
G\\\text{is bijective}}}b_{f\left(  1\right)  ,1}b_{f\left(  2\right)
,2}\cdots b_{f\left(  n\right)  ,n}.
\]

\end{exercise}

\begin{remark}
Let $\mathbb{K}$ be a commutative ring. Let $n\in\mathbb{N}$. Let $A=\left(
a_{i,j}\right)  _{1\leq i\leq n,\ 1\leq j\leq n}\in\mathbb{K}^{n\times n}$ be
an $n\times n$-matrix. Then, the \textit{permanent} $\operatorname*{per}A$ of
$A$ is defined to be the element%
\[
\sum_{\sigma\in S_{n}}a_{1,\sigma\left(  1\right)  }a_{2,\sigma\left(
2\right)  }\cdots a_{n,\sigma\left(  n\right)  }%
\]
of $\mathbb{K}$. The concept of a permanent is thus similar to the concept of
a determinant (the only difference is that the factor $\left(  -1\right)
^{\sigma}$ is missing from the definition of the permanent); however, it has
far fewer interesting properties. One of the properties that it does have is
the so-called \textit{Ryser formula} (see, e.g., \cite[{\S 4.9, [9e]}%
]{Comtet74}), which says that%
\[
\operatorname*{per}A=\left(  -1\right)  ^{n}\sum_{I\subseteq\left\{
1,2,\ldots,n\right\}  }\left(  -1\right)  ^{\left\vert I\right\vert }%
\prod_{j=1}^{n}\sum_{i\in I}a_{i,j}.
\]
The reader is invited to check that this formula follows from Exercise
\ref{exe.det.sumdets1} \textbf{(c)} (applied to $\mathbb{L}=\mathbb{K}$,
$G=\left\{  1,2,\ldots,n\right\}  $, $H=\varnothing$ and $b_{i,j}=a_{i,j}$).
\end{remark}

\begin{exercise}
\label{exe.det.sumdets2}Let $n\in\mathbb{N}$. Let $G$ be a finite set such
that $n<\left\vert G\right\vert $. For each $i\in G$, let $A_{i}\in
\mathbb{K}^{n\times n}$ be an $n\times n$-matrix. Prove that%
\[
\sum_{I\subseteq G}\left(  -1\right)  ^{\left\vert I\right\vert }\det\left(
\sum_{i\in I}A_{i}\right)  =0.
\]

\end{exercise}

\begin{exercise}
\label{exe.powerdet.gen}Let $n\in\mathbb{N}$. Let $A=\left(  a_{i,j}\right)
_{1\leq i\leq n,\ 1\leq j\leq n}$ be an $n\times n$-matrix.

\textbf{(a)} Prove that%
\[
\sum_{\sigma\in S_{n}}\left(  -1\right)  ^{\sigma}\left(  \sum_{i=1}%
^{n}a_{i,\sigma\left(  i\right)  }\right)  ^{k}=0
\]
for each $k\in\left\{  0,1,\ldots,n-2\right\}  $.

\textbf{(b)} Assume that $n\geq1$. Prove that%
\[
\sum_{\sigma\in S_{n}}\left(  -1\right)  ^{\sigma}\left(  \sum_{i=1}%
^{n}a_{i,\sigma\left(  i\right)  }\right)  ^{n-1}=\left(  n-1\right)
!\cdot\sum_{p=1}^{n}\sum_{q=1}^{n}\left(  -1\right)  ^{p+q}\det\left(  A_{\sim
p,\sim q}\right)  .
\]
(Here, we are using the notation introduced in Definition
\ref{def.submatrix.minor}.)
\end{exercise}

\begin{exercise}
\label{exe.powerdet.r1}Let $n\in\mathbb{N}$. Let $a_{1},a_{2},\ldots,a_{n}$ be
$n$ elements of $\mathbb{K}$. Let $b_{1},b_{2},\ldots,b_{n}$ be $n$ elements
of $\mathbb{K}$.

Let $m=\dbinom{n}{2}$.

\textbf{(a)} Prove that%
\[
\sum_{\sigma\in S_{n}}\left(  -1\right)  ^{\sigma}\left(  \sum_{i=1}^{n}%
a_{i}b_{\sigma\left(  i\right)  }\right)  ^{k}=0
\]
for each $k\in\left\{  0,1,\ldots,m-1\right\}  $.

\textbf{(b)} Prove that%
\[
\sum_{\sigma\in S_{n}}\left(  -1\right)  ^{\sigma}\left(  \sum_{i=1}^{n}%
a_{i}b_{\sigma\left(  i\right)  }\right)  ^{m}=\mathbf{m}\left(
0,1,\ldots,n-1\right)  \cdot\prod_{1\leq i<j\leq n}\left(  \left(  a_{i}%
-a_{j}\right)  \left(  b_{i}-b_{j}\right)  \right)  .
\]
Here, we are using the notation introduced in Exercise \ref{exe.multinom2}.

\textbf{(c)} Let $k\in\mathbb{N}$. Prove that%
\begin{align*}
&  \sum_{\sigma\in S_{n}}\left(  -1\right)  ^{\sigma}\left(  \sum_{i=1}%
^{n}a_{i}b_{\sigma\left(  i\right)  }\right)  ^{k}\\
&  =\sum_{\substack{\left(  g_{1},g_{2},\ldots,g_{n}\right)  \in\mathbb{N}%
^{n};\\g_{1}<g_{2}<\cdots<g_{n};\\g_{1}+g_{2}+\cdots+g_{n}=k}}\mathbf{m}%
\left(  g_{1},g_{2},\ldots,g_{n}\right) \\
&  \ \ \ \ \ \ \ \ \ \ \cdot\det\left(  \left(  a_{i}^{g_{j}}\right)  _{1\leq
i\leq n,\ 1\leq j\leq n}\right)  \cdot\det\left(  \left(  b_{i}^{g_{j}%
}\right)  _{1\leq i\leq n,\ 1\leq j\leq n}\right)  .
\end{align*}
Here, we are using the notation introduced in Exercise \ref{exe.multinom2}.
\end{exercise}

Note that Exercise \ref{exe.powerdet.r1} generalizes \cite[Exercise
12.13]{AndDos} and \cite[\S 12.1, Problem 1]{AndDosS}.

\subsection{Additional exercises}

Here are a few more additional exercises, with no importance to the rest of
the text (and mostly no solutions given).\Needspace{8cm}

\begin{exercise}
\label{addexe.jacobi-complement}Let $n\in\mathbb{N}$. For any subset $I$ of
$\left\{  1,2,\ldots,n\right\}  $, we let $\widetilde{I}$ denote the
complement $\left\{  1,2,\ldots,n\right\}  \setminus I$ of $I$. (For instance,
if $n=4$ and $I=\left\{  1,4\right\}  $, then $\widetilde{I}=\left\{
2,3\right\}  $.) Let us use the notations introduced in Definition
\ref{def.submatrix} and in Definition \ref{def.sect.laplace.notations}.

Let $A\in\mathbb{K}^{n\times n}$ be an invertible matrix. Let $P$ and $Q$ be
two subsets of $\left\{  1,2,\ldots,n\right\}  $ satisfying $\left\vert
P\right\vert =\left\vert Q\right\vert $. Prove that%
\begin{equation}
\det\left(  \operatorname*{sub}\nolimits_{w\left(  P\right)  }^{w\left(
Q\right)  }A\right)  =\left(  -1\right)  ^{\sum P+\sum Q}\det A\cdot
\det\left(  \operatorname*{sub}\nolimits_{w\left(  \widetilde{Q}\right)
}^{w\left(  \widetilde{P}\right)  }\left(  A^{-1}\right)  \right)  .
\label{eq.addexe.jacobi-complement.1}%
\end{equation}

[\textbf{Hint:} Apply Exercise \ref{exe.block2x2.jacobi.rewr} to a matrix
obtained from $A$ by permuting the rows and permuting the columns.]
\end{exercise}

\begin{remark}
The claim of Exercise \ref{addexe.jacobi-complement} is the so-called
\textit{Jacobi complementary minor theorem}. It appears, for example, in
\cite[(1)]{Lalonde} and in \cite[Lemma A.1 (e)]{CaSoSp12}, and is used rather
often when working with determinants (for example, it is used in \cite[Chapter
SYM, proof of Proposition (7.5) (5)]{LLPT95} and many times in \cite{CaSoSp12}).

The determinant of a submatrix of a matrix $A$ is called a \textit{minor} of
$A$. Thus, in the equality (\ref{eq.addexe.jacobi-complement.1}), the
determinant $\det\left(  \operatorname*{sub}\nolimits_{w\left(  P\right)
}^{w\left(  Q\right)  }A\right)  $ on the left hand side is a minor of $A$,
whereas the determinant $\det\left(  \operatorname*{sub}\nolimits_{w\left(
\widetilde{Q}\right)  }^{w\left(  \widetilde{P}\right)  }\left(
A^{-1}\right)  \right)  $ on the right hand side is a minor of $A^{-1}$. Thus,
roughly speaking, the equality (\ref{eq.addexe.jacobi-complement.1}) says that
any minor of $A$ equals a certain minor of $A^{-1}$ times $\det A$ times a
certain sign.

It is instructive to check the particular case of
(\ref{eq.addexe.jacobi-complement.1}) obtained when both $P$ and $Q$ are sets
of cardinality $n-1$ (so that $\widetilde{P}$ and $\widetilde{Q}$ are
$1$-element sets). This particular case turns out to be the statement of
Theorem \ref{thm.matrices.inverses.square} \textbf{(b)} in disguise.

Exercise \ref{exe.block2x2.jacobi.rewr} is the particular case of Exercise
\ref{addexe.jacobi-complement} obtained when $P=\left\{  1,2,\ldots,k\right\}
$ and $Q=\left\{  1,2,\ldots,k\right\}  $.
\end{remark}

\begin{exercise}
\label{exe.det.pluecker.multi}Let $n$ and $k$ be positive integers such that
$k\leq n$. Let $A\in\mathbb{K}^{n\times\left(  n-k\right)  }$ and
$B\in\mathbb{K}^{n\times\left(  n+k\right)  }$.

Let us use the notations from Definition \ref{def.unrows}. For any subset $I$
of $\left\{  1,2,\ldots,n+k\right\}  $, we introduce the following five notations:

\begin{itemize}
\item Let $\sum I$ denote the sum of all elements of $I$. (Thus, $\sum
I=\sum_{i\in I}i$.)

\item Let $w\left(  I\right)  $ denote the list of all elements of $I$ in
increasing order (with no repetitions). (See Definition \ref{def.ind.inclist}
for the formal definition of this list.) (For example, $w\left(  \left\{
3,4,8\right\}  \right)  =\left(  3,4,8\right)  $.)

\item Let $\left(  A\mid B_{\bullet,I}\right)  $ denote the $n\times\left(
n-k+\left\vert I\right\vert \right)  $-matrix whose columns are
$\underbrace{A_{\bullet,1},A_{\bullet,2},\ldots,A_{\bullet,n-k}}_{\text{the
columns of }A},B_{\bullet,i_{1}},B_{\bullet,i_{2}},\ldots,B_{\bullet,i_{\ell}%
}$ (from left to right), where $\left(  i_{1},i_{2},\ldots,i_{\ell}\right)
=w\left(  I\right)  $.

\item Let $B_{\bullet,\sim I}$ denote the $n\times\left(  n+k-\left\vert
I\right\vert \right)  $-matrix whose columns are $B_{\bullet,j_{1}}%
,B_{\bullet,j_{2}},\ldots,B_{\bullet,j_{h}}$ (from left to right), where
$\left(  j_{1},j_{2},\ldots,j_{h}\right)  =w\left(  \left\{  1,2,\ldots
,n+k\right\}  \setminus I\right)  $. \newline(Using the notations of
Definition \ref{def.rowscols}, we can rewrite this definition as
$B_{\bullet,\sim I}=\operatorname*{cols}\nolimits_{j_{1},j_{2},\ldots,j_{h}}%
B$, where $\left(  j_{1},j_{2},\ldots,j_{h}\right)  =w\left(  \left\{
1,2,\ldots,n+k\right\}  \setminus I\right)  $.)
\end{itemize}

Then, prove that%
\[
\sum_{\substack{I\subseteq\left\{  1,2,\ldots,n+k\right\}  ;\\\left\vert
I\right\vert =k}}\left(  -1\right)  ^{\sum I+\left(  1+2+\cdots+k\right)
}\det\left(  A\mid B_{\bullet,I}\right)  \det\left(  B_{\bullet,\sim
I}\right)  =0.
\]
(Note that this generalizes Theorem \ref{thm.pluecker.plu}; indeed, the latter
theorem is the particular case for $k=1$.)
\end{exercise}

\begin{exercise}
\label{exe.det.delannoy}Recall that the binomial coefficients satisfy the
recurrence relation (\ref{eq.binom.rec.m}), which (visually) says that every
entry of Pascal's triangle is the sum of the two entries left-above it and
right-above it.

Let us now define a variation of Pascal's triangle as follows: Define a
nonnegative integer $\dbinom{m}{n}_{D}$ for every $m\in\mathbb{N}$ and
$n\in\mathbb{N}$ recursively as follows:

\begin{itemize}
\item Set $\dbinom{0}{n}_{D}=%
\begin{cases}
1, & \text{if }n=0;\\
0, & \text{if }n>0
\end{cases}
$ for every $n\in\mathbb{N}$.

\item For every $m\in\mathbb{Z}$ and $n\in\mathbb{Z}$, set $\dbinom{m}{n}%
_{D}=0$ if either $m$ or $n$ is negative.

\item For every positive integer $m$ and every $n\in\mathbb{N}$, set
\[
\dbinom{m}{n}_{D}=\dbinom{m-1}{n-1}_{D}+\dbinom{m-1}{n}_{D}+\dbinom{m-2}%
{n-1}_{D}.
\]

\end{itemize}

\noindent(Thus, if we lay these $\dbinom{m}{n}_{D}$ out in the same way as the
binomial coefficients $\dbinom{m}{n}$ in Pascal's triangle, then every entry
is the sum of the three entries left-above it, right-above it, and straight
above it.)

The integers $\dbinom{m}{n}_{D}$ are known as the
\href{https://en.wikipedia.org/wiki/Delannoy_number}{Delannoy numbers}.

\textbf{(a)} Show that
\[
\dbinom{n+m}{n}_{D}=\sum_{i=0}^{n}\dbinom{n}{i}\dbinom{m+i}{n}=\sum_{i=0}%
^{n}\dbinom{n}{i}\dbinom{m}{i}2^{i}%
\]
for every $n\in\mathbb{N}$ and $m\in\mathbb{N}$. (The second equality sign
here is a consequence of Proposition \ref{prop.binom.bin-id} \textbf{(e)}.)

\textbf{(b)} Let $n\in\mathbb{N}$. Let $A$ be the $n\times n$-matrix $\left(
\dbinom{i+j-2}{i-1}_{D}\right)  _{1\leq i\leq n,\ 1\leq j\leq n}$ (an analogue
of the matrix $A$ from Exercise \ref{exe.ps4.pascal}). Show that%
\[
\det A=2^{n\left(  n-1\right)  /2}.
\]

\end{exercise}

\begin{noncompile}
\textbf{(a)} One can prove $\dbinom{n+m}{n}_{D}=\sum_{i=0}^{n}\dbinom{n}%
{i}\dbinom{m+i}{n}$ by strong induction over $n+m$, along the following lines:%
\begin{align*}
&  \sum_{i=0}^{n}\dbinom{n}{i}\underbrace{\dbinom{m+i}{n}}_{=\dbinom
{m+i-1}{n-1}+\dbinom{m+i-1}{n}}\\
&  =\sum_{i=0}^{n}\underbrace{\dbinom{n}{i}}_{=\dbinom{n-1}{i-1}+\dbinom
{n-1}{i}}\dbinom{m+i-1}{n-1}+\underbrace{\sum_{i=0}^{n}\dbinom{n}{i}%
\dbinom{m+i-1}{n}}_{=\dbinom{n+m-1}{n}_{D}}\\
&  =\underbrace{\sum_{i=0}^{n}\dbinom{n-1}{i-1}\dbinom{m+i-1}{n-1}%
}_{\substack{=\sum_{i=0}^{n}\dbinom{n-1}{i-1}\dbinom{m+\left(  i-1\right)
}{n-1}\\=\dbinom{\left(  n-1\right)  +m}{n-1}_{D}}}+\underbrace{\sum_{i=0}%
^{n}\dbinom{n-1}{i}\dbinom{m+i-1}{n-1}}_{\substack{=\sum_{i=0}^{n}\dbinom
{n-1}{i}\dbinom{\left(  m-1\right)  +i}{n-1}\\=\dbinom{\left(  n-1\right)
+\left(  m-1\right)  }{n-1}_{D}}}+\dbinom{n+m-1}{n}_{D}\\
&  =\dbinom{\left(  n-1\right)  +m}{n-1}_{D}+\dbinom{\left(  n-1\right)
+\left(  m-1\right)  }{n-1}_{D}+\dbinom{n+m-1}{n}_{D}=\dbinom{n+m}{n}_{D}.
\end{align*}

\textbf{(b)} Use $\dbinom{n+m}{n}_{D}=\sum_{i=0}^{n}\dbinom{n}{i}\dbinom{m}%
{i}2^{i}$ and the same tactic as in Exercise \ref{exe.ps4.pascal}.
\end{noncompile}

The next exercise is known under the (arguably not very distinctive) name of
\textquotedblleft%
\href{https://en.wikipedia.org/wiki/Matrix_determinant_lemma}{matrix
determinant lemma}\textquotedblright:

\begin{exercise}
\label{exe.det.rk1upd}Let $n\in\mathbb{N}$. Let $u$ be a column vector with
$n$ entries, and let $v$ be a row vector with $n$ entries. (Thus, $uv$ is an
$n\times n$-matrix, whereas $vu$ is a $1\times1$-matrix.) Let $A$ be an
$n\times n$-matrix. Prove that%
\[
\det\left(  A+uv\right)  =\det A+v\left(  \operatorname*{adj}A\right)  u
\]
(where we regard the $1\times1$-matrix $v\left(  \operatorname*{adj}A\right)
u$ as an element of $\mathbb{K}$).
\end{exercise}

The next exercise relies on Definition \ref{def.block2x2}.

\begin{exercise}
\label{exe.det.bordered}Let $n\in\mathbb{N}$. Let $u\in\mathbb{K}^{n\times1}$
be a column vector with $n$ entries, and let $v\in\mathbb{K}^{1\times n}$ be a
row vector with $n$ entries. (Thus, $uv$ is an $n\times n$-matrix, whereas
$vu$ is a $1\times1$-matrix.) Let $h\in\mathbb{K}$. Let $H$ be the $1\times
1$-matrix $\left(
\begin{array}
[c]{c}%
h
\end{array}
\right)  \in\mathbb{K}^{1\times1}$.

\textbf{(a)} Prove that every $n\times n$-matrix $A\in\mathbb{K}^{n\times n}$
satisfies%
\[
\det\left(
\begin{array}
[c]{cc}%
A & u\\
v & H
\end{array}
\right)  =h\det A-v\left(  \operatorname*{adj}A\right)  u
\]
(where we regard the $1\times1$-matrix $v\left(  \operatorname*{adj}A\right)
u$ as an element of $\mathbb{K}$).

\textbf{(b)} Write the vector $u$ in the form $u=\left(  u_{1},u_{2}%
,\ldots,u_{n}\right)  ^{T}$.\ Write the vector $v$ in the form $v=\left(
v_{1},v_{2},\ldots,v_{n}\right)  $.

For every two objects $i$ and $j$, define $\delta_{i,j}\in\mathbb{K}$ by
$\delta_{i,j}=%
\begin{cases}
1, & \text{if }i=j;\\
0, & \text{if }i\neq j
\end{cases}
$.

Let $d_{1},d_{2},\ldots,d_{n}$ be $n$ elements of $\mathbb{K}$. Let $D$ be the
$n\times n$-matrix $\left(  d_{i}\delta_{i,j}\right)  _{1\leq i\leq n,\ 1\leq
j\leq n}$. Prove that%
\[
\det\left(
\begin{array}
[c]{cc}%
D & u\\
v & H
\end{array}
\right)  =h\cdot\left(  d_{1}d_{2}\cdots d_{n}\right)  -\sum_{i=1}^{n}%
u_{i}v_{i}\prod_{\substack{j\in\left\{  1,2,\ldots,n\right\}  ;\\j\neq
i}}d_{j}.
\]

\end{exercise}

\begin{example}
Let $n=3$. Then, Exercise \ref{exe.det.bordered} \textbf{(b)} states that
\begin{align*}
\det\left(
\begin{array}
[c]{cccc}%
d_{1} & 0 & 0 & u_{1}\\
0 & d_{2} & 0 & u_{2}\\
0 & 0 & d_{3} & u_{3}\\
v_{1} & v_{2} & v_{3} & h
\end{array}
\right)   &  =h\cdot\left(  d_{1}d_{2}d_{3}\right)  -\sum_{i=1}^{3}u_{i}%
v_{i}\prod_{\substack{j\in\left\{  1,2,3\right\}  ;\\j\neq i}}d_{j}\\
&  =hd_{1}d_{2}d_{3}-\left(  u_{1}v_{1}d_{2}d_{3}+u_{2}v_{2}d_{1}d_{3}%
+u_{3}v_{3}d_{1}d_{2}\right)
\end{align*}
for any ten elements $u_{1},u_{2},u_{3},v_{1},v_{2},v_{3},d_{1},d_{2},d_{3},h$
of $\mathbb{K}$.
\end{example}

\begin{exercise}
\label{exe.det.resultant}Let $P=\sum_{k=0}^{d}p_{k}X^{k}$ and $Q=\sum
_{k=0}^{e}q_{k}X^{k}$ be two polynomials over $\mathbb{K}$ (where $p_{0}%
,p_{1},\ldots,p_{d}\in\mathbb{K}$ and $q_{0},q_{1},\ldots,q_{e}\in\mathbb{K}$
are their coefficients) such that $d+e>0$. Define a $\left(  d+e\right)
\times\left(  d+e\right)  $-matrix $A$ as follows:

\begin{itemize}
\item For every $k\in\left\{  1,2,\ldots,e\right\}  $, the $k$-th row of $A$
is%
\[
\left(  \underbrace{0,0,\ldots,0}_{k-1\text{ zeroes}},p_{d},p_{d-1}%
,\ldots,p_{1},p_{0},\underbrace{0,0,\ldots,0}_{e-k\text{ zeroes}}\right)  .
\]

\item For every $k\in\left\{  1,2,\ldots,d\right\}  $, the $\left(
e+k\right)  $-th row of $A$ is%
\[
\left(  \underbrace{0,0,\ldots,0}_{k-1\text{ zeroes}},q_{e},q_{e-1}%
,\ldots,q_{1},q_{0},\underbrace{0,0,\ldots,0}_{d-k\text{ zeroes}}\right)  .
\]

\end{itemize}

(For example, if $d=4$ and $e=3$, then%
\[
A=\left(
\begin{array}
[c]{ccccccc}%
p_{4} & p_{3} & p_{2} & p_{1} & p_{0} & 0 & 0\\
0 & p_{4} & p_{3} & p_{2} & p_{1} & p_{0} & 0\\
0 & 0 & p_{4} & p_{3} & p_{2} & p_{1} & p_{0}\\
q_{3} & q_{2} & q_{1} & q_{0} & 0 & 0 & 0\\
0 & q_{3} & q_{2} & q_{1} & q_{0} & 0 & 0\\
0 & 0 & q_{3} & q_{2} & q_{1} & q_{0} & 0\\
0 & 0 & 0 & q_{3} & q_{2} & q_{1} & q_{0}%
\end{array}
\right)  .
\]
)

Assume that the polynomials $P$ and $Q$ have a common root $z$ (that is, there
exists a $z\in\mathbb{K}$ such that $P\left(  z\right)  =0$ and $Q\left(
z\right)  =0$). Show that $\det A=0$.

[\textbf{Hint:} Find a column vector $v$ with $d+e$ entries satisfying
$Av=0_{\left(  d+e\right)  \times1}$; then apply Corollary
\ref{cor.adj.kernel}.]
\end{exercise}

\begin{remark}
The matrix $A$ in Exercise \ref{exe.det.resultant} is called the
\href{https://en.wikipedia.org/wiki/Sylvester_matrix}{\textit{Sylvester
matrix}} of the polynomials $P$ and $Q$ (for degrees $d$ and $e$); its
determinant $\det A$ is known as their
\href{https://en.wikipedia.org/wiki/Resultant}{\textit{resultant}} (at least
when $d$ and $e$ are actually the degrees of $P$ and $Q$). According to the
exercise, the condition $\det A=0$ is necessary for $P$ and $Q$ to have a
common root. In the general case, the converse does not hold: For one, you can
always force $\det A$ to be $0$ by taking $d>\deg P$ and $e>\deg Q$ (so
$p_{d}=0$ and $q_{e}=0$, and thus the $1$-st column of $A$ consists of
zeroes). More importantly, the resultant of the two polynomials $X^{3}-1$ and
$X^{2}+X+1$ is $0$, but they only have common roots in $\mathbb{C}$, not in
$\mathbb{R}$. Thus, there is more to common roots than just the vanishing of a determinant.

However, if $\mathbb{K}$ is an algebraically closed field (I won't go into the
details of what this means, but an example of such a field is $\mathbb{C}$),
and if $d=\deg P$ and $e=\deg Q$, then the polynomials $P$ and $Q$ have a
common root \textbf{if and only if} their resultant is $0$.
\end{remark}

The next two exercises can be seen as variations on the Vandermonde
determinant. The first one is \cite[Proposition 1]{Krattenthaler} (and also
appears in \cite[Theorem 2]{GriHyp}):

\begin{exercise}
\label{exe.det.vdm-pol}Let $n\in\mathbb{N}$. Let $a_{1},a_{2},\ldots,a_{n}%
\in\mathbb{K}$. Furthermore, for each $j\in\left\{  1,2,\ldots,n\right\}  $,
let $P_{j}\in\mathbb{K}\left[  X\right]  $ be a polynomial such that%
\[
\deg\left(  P_{j}\right)  \leq j-1.
\]
(In particular, $\deg\left(  P_{1}\right)  \leq1-1=0$, so that the polynomial
$P_{1}$ is constant.)

For each $j\in\left\{  1,2,\ldots,n\right\}  $, let $c_{j}$ be the coefficient
of $X^{j-1}$ in the polynomial $P_{j}$. Prove that
\[
\det\left(  \left(  P_{j}\left(  a_{i}\right)  \right)  _{1\leq i\leq
n,\ 1\leq j\leq n}\right)  =\left(  \prod_{j=1}^{n}c_{j}\right)  \cdot
\prod_{1\leq j<i\leq n}\left(  a_{i}-a_{j}\right)  .
\]

\end{exercise}

The next exercise is \cite[Lemma 6]{Krattenthaler} (with the notations changed):

\begin{exercise}
\label{exe.det.kratt-lem6}Let $n\in\mathbb{N}$. Let $a_{1},a_{2},\ldots
,a_{n}\in\mathbb{K}$ and $b_{2},b_{3},\ldots,b_{n}\in\mathbb{K}$. For each
$j\in\left\{  1,2,\ldots,n\right\}  $, we define a polynomial $Q_{j}%
\in\mathbb{K}\left[  X\right]  $ by setting%
\[
Q_{j}=\left(  X-b_{j+1}\right)  \left(  X-b_{j+2}\right)  \cdots\left(
X-b_{n}\right)  .
\]
(In particular, $Q_{n}=\left(  \text{empty product}\right)  =1$.) Furthermore,
for each $j\in\left\{  1,2,\ldots,n\right\}  $, let $P_{j}\in\mathbb{K}\left[
X\right]  $ be a polynomial such that%
\[
\deg\left(  P_{j}\right)  \leq j-1.
\]
(In particular, $\deg\left(  P_{1}\right)  \leq1-1=0$, so that the polynomial
$P_{1}$ is constant.)

Prove that
\[
\det\left(  \left(  P_{j}\left(  a_{i}\right)  Q_{j}\left(  a_{i}\right)
\right)  _{1\leq i\leq n,\ 1\leq j\leq n}\right)  =\left(  \prod_{j=1}%
^{n}P_{j}\left(  b_{j}\right)  \right)  \cdot\prod_{1\leq i<j\leq n}\left(
a_{i}-a_{j}\right)  .
\]

[\textbf{Hint:} For each $j\in\left\{  1,2,\ldots,n\right\}  $, construct $j$
elements $c_{j,1},c_{j,2},\ldots,c_{j,j}$ of $\mathbb{K}$ satisfying
$P_{j}Q_{j}=\sum_{k=1}^{j}c_{j,k}Q_{k}$ and $c_{j,j}=P_{j}\left(
b_{j}\right)  $.]
\end{exercise}

Exercise \ref{exe.det.kratt-lem6} has many applications:\footnote{Exercise
\ref{exe.det.kratt-lem6-app} \textbf{(b)} appears in \cite[Lemma
3]{Krattenthaler}, in \cite[Lemma 2.7]{Gorin21}, in \cite[\S 1.2.3, exercise
47]{Knuth-TAoCP1} and in \cite[Theorem 2.9]{Bresso99}. (In all of these
sources, it appears with modified notations. For instance, \cite[Lemma
3]{Krattenthaler} is Exercise \ref{exe.det.kratt-lem6-app} \textbf{(b)} for
$a_{i}=X_{i}$, $b_{i}=A_{i}$ and $c_{i}=B_{i+1}$.)}

\begin{exercise}
\label{exe.det.kratt-lem6-app}Let $n\in\mathbb{N}$. Let $a_{1},a_{2}%
,\ldots,a_{n}\in\mathbb{K}$ and $b_{2},b_{3},\ldots,b_{n}\in\mathbb{K}$ and
$c_{1},c_{2},\ldots,c_{n-1}\in\mathbb{K}$.

\textbf{(a)} Prove that
\begin{align*}
&  \det\left(  \left(  \left(  \prod_{u=1}^{j-1}\left(  a_{i}-c_{u}\right)
\right)  \left(  \prod_{u=j+1}^{n}\left(  a_{i}-b_{u}\right)  \right)
\right)  _{1\leq i\leq n,\ 1\leq j\leq n}\right) \\
&  =\left(  \prod_{1\leq i<j\leq n}\left(  b_{j}-c_{i}\right)  \right)
\cdot\prod_{1\leq i<j\leq n}\left(  a_{i}-a_{j}\right)  .
\end{align*}

\textbf{(b)} Prove that
\begin{align*}
&  \det\left(  \left(  \left(  \prod_{u=1}^{j-1}\left(  a_{i}+c_{u}\right)
\right)  \left(  \prod_{u=j+1}^{n}\left(  a_{i}+b_{u}\right)  \right)
\right)  _{1\leq i\leq n,\ 1\leq j\leq n}\right) \\
&  =\left(  \prod_{1\leq i<j\leq n}\left(  c_{i}-b_{j}\right)  \right)
\cdot\prod_{1\leq i<j\leq n}\left(  a_{i}-a_{j}\right)  .
\end{align*}

\textbf{(c)} Use this to solve Exercise \ref{exe.cauchy-det} again.
\end{exercise}

The next exercise (more precisely, its part \textbf{(a)}) is a recent result
of Fraser and Yeats \cite[Theorem 2.1.4]{FraYea21} (slightly extended):

\begin{exercise}
\label{exe.det.fraser-yeats}\textbf{(a)} Let $n\in\mathbb{N}$. Let
$A\in\mathbb{K}^{n\times n}$ be an $n\times n$-matrix. Let $k\in\mathbb{N}$.
Let $u_{1},u_{2},\ldots,u_{k}$ be $k$ elements of $\left\{  1,2,\ldots
,n\right\}  $. Let $v_{1},v_{2},\ldots,v_{k}$ be $k$ further elements of
$\left\{  1,2,\ldots,n\right\}  $. Let $U=\left\{  u_{1},u_{2},\ldots
,u_{k}\right\}  $ and $V=\left\{  v_{1},v_{2},\ldots,v_{k}\right\}  $. Let
$t\in\left\{  1,2,\ldots,n\right\}  \setminus V$. Prove that%
\begin{align*}
&  \det A\cdot\det\left(  \operatorname*{sub}\nolimits_{u_{1},u_{2}%
,\ldots,u_{k}}^{v_{1},v_{2},\ldots,v_{k}}A\right) \\
&  =\sum_{s\in\left\{  1,2,\ldots,n\right\}  \setminus U}\left(  -1\right)
^{s+t}\det\left(  A_{\sim s,\sim t}\right)  \det\left(  \operatorname*{sub}%
\nolimits_{s,u_{1},u_{2},\ldots,u_{k}}^{t,v_{1},v_{2},\ldots,v_{k}}A\right)  .
\end{align*}
(Here, we are using the notations of Definition \ref{def.submatrix} and of
Definition \ref{def.submatrix.minor}.)

\textbf{(b)} Use this to give a new proof of Proposition
\ref{prop.desnanot.12}.
\end{exercise}

For the next two exercises, we need a notation for exchanging certain columns
between two matrices:

\begin{definition}
\label{def.matrix-col-exchange}Let $n\in\mathbb{N}$ and $m\in\mathbb{N}$ and
$p\in\mathbb{N}$. Let $A\in\mathbb{K}^{n\times m}$ and $B\in\mathbb{K}%
^{n\times p}$ be two matrices. Let $J$ be a subset of $\left\{  1,2,\ldots
,m\right\}  $, and let $K$ be a subset of $\left\{  1,2,\ldots,p\right\}  $
such that $\left\vert J\right\vert =\left\vert K\right\vert $. Let
$j_{1},j_{2},\ldots,j_{q}$ be all elements of $J$, listed in increasing order
(with no repetitions). Thus, $q=\left\vert J\right\vert =\left\vert
K\right\vert $. Hence, the set $K$ has exactly $q$ elements. Let $k_{1}%
,k_{2},\ldots,k_{q}$ be all elements of $K$, listed in increasing order (with
no repetitions). (This is well-defined, since $K$ has exactly $q$ elements.)
Then, we let $\left(
\begin{array}
[c]{c}%
A\leftarrow B\\
J\leftarrow K
\end{array}
\right)  $ denote the $n\times m$-matrix obtained from $A$ by replacing the
$j_{1}$-st, $j_{2}$-nd, $\ldots$, $j_{q}$-th columns of $A$ by the $k_{1}$-st,
$k_{2}$-nd, $\ldots$, $k_{q}$-th columns of $B$, respectively. In other words,
if we write the matrices $A$ and $B$ as $A=\left(  a_{x,y}\right)  _{1\leq
x\leq n,\ 1\leq y\leq m}$ and $B=\left(  b_{x,y}\right)  _{1\leq x\leq
n,\ 1\leq y\leq p}$, then $\left(
\begin{array}
[c]{c}%
A\leftarrow B\\
J\leftarrow K
\end{array}
\right)  $ is defined to be the $n\times m$-matrix $\left(  c_{x,y}\right)
_{1\leq x\leq n,\ 1\leq y\leq m}$, where we let%
\begin{align}
c_{x,y}  &  =%
\begin{cases}
a_{x,y}, & \text{if }y\notin J;\\
b_{x,k_{i}}, & \text{if }y=j_{i}\text{ for some }i\in\left\{  1,2,\ldots
,q\right\}
\end{cases}
\label{eq.def.matrix-col-exchange.cxy=}\\
&  \ \ \ \ \ \ \ \ \ \ \ \ \ \ \ \ \ \ \ \ \text{for all }x\in\left\{
1,2,\ldots,n\right\}  \text{ and }y\in\left\{  1,2,\ldots,m\right\}
.\nonumber
\end{align}

\end{definition}

\begin{example}
Let $A=\left(
\begin{array}
[c]{cccc}%
a & b & c & d\\
a^{\prime} & b^{\prime} & c^{\prime} & d^{\prime}\\
a^{\prime\prime} & b^{\prime\prime} & c^{\prime\prime} & d^{\prime\prime}%
\end{array}
\right)  $ and $B=\left(
\begin{array}
[c]{ccc}%
e & f & g\\
e^{\prime} & f^{\prime} & g^{\prime}\\
e^{\prime\prime} & f^{\prime\prime} & g^{\prime\prime}%
\end{array}
\right)  $ and $J=\left\{  1,3\right\}  $ and $K=\left\{  2,3\right\}  $.
Then,%
\begin{align*}
\left(
\begin{array}
[c]{c}%
A\leftarrow B\\
J\leftarrow K
\end{array}
\right)   &  =\left(
\begin{array}
[c]{cccc}%
f & b & g & d\\
f^{\prime} & b^{\prime} & g^{\prime} & d^{\prime}\\
f^{\prime\prime} & b^{\prime\prime} & g^{\prime\prime} & d^{\prime\prime}%
\end{array}
\right)  \ \ \ \ \ \ \ \ \ \ \text{and}\\
\left(
\begin{array}
[c]{c}%
B\leftarrow A\\
K\leftarrow J
\end{array}
\right)   &  =\left(
\begin{array}
[c]{ccc}%
e & a & c\\
e^{\prime} & a^{\prime} & c^{\prime}\\
e^{\prime\prime} & a^{\prime\prime} & c^{\prime\prime}%
\end{array}
\right)  .
\end{align*}

\end{example}

The next exercise is known as \textit{Sylvester's Lemma} (appearing, e.g., in
\cite[Chapter 8, Lemma 2]{Fulton-Young} or \cite[\S 137]{Muir}):

\begin{exercise}
\label{exe.det.syl-lem}Let $n\in\mathbb{N}$. Let $A\in\mathbb{K}^{n\times n}$
and $B\in\mathbb{K}^{n\times n}$ be two $n\times n$-matrices, and let $J$ be a
subset of $\left\{  1,2,\ldots,n\right\}  $. Prove that
\[
\det A\cdot\det B=\sum_{\substack{K\subseteq\left\{  1,2,\ldots,n\right\}
;\\\left\vert K\right\vert =\left\vert J\right\vert }}\det\left(
\begin{array}
[c]{c}%
A\leftarrow B\\
J\leftarrow K
\end{array}
\right)  \det\left(
\begin{array}
[c]{c}%
B\leftarrow A\\
K\leftarrow J
\end{array}
\right)  .
\]

\end{exercise}

The next exercise is essentially \cite[\S 319]{Muir}:

\begin{exercise}
\label{exe.det.muir-col-by-col}Let $n\in\mathbb{N}$. Let $A\in\mathbb{K}%
^{n\times n}$ and $B\in\mathbb{K}^{n\times n}$ be two $n\times n$-matrices,
and let $r\in\mathbb{N}$. Prove that
\[
\sum_{\substack{J\subseteq\left\{  1,2,\ldots,n\right\}  ;\\\left\vert
J\right\vert =r}}\det\left(
\begin{array}
[c]{c}%
A\leftarrow B\\
J\leftarrow J
\end{array}
\right)  =\sum_{\substack{J\subseteq\left\{  1,2,\ldots,n\right\}
;\\\left\vert J\right\vert =r}}\det\left(
\begin{array}
[c]{c}%
A^{T}\leftarrow B^{T}\\
J\leftarrow J
\end{array}
\right)  .
\]

\end{exercise}

The next exercise generalizes \cite[\S 231]{Muir}\footnote{Muir, in
\cite[\S 231]{Muir}, requires $\mathbb{K}$ to be the field $\mathbb{R}$ of
real numbers, and (rather surprisingly) his proof actually uses this
(extraneous) requirement.} as well as \cite[\S VI.4, ASSERTION]{AlAmra97}%
\footnote{To wit, \cite[\S VI.4, ASSERTION]{AlAmra97} can be obtained by
applying Exercise \ref{exe.det.AB=0-proport} to $k=1$, $A=\left(
U_{i,j}\right)  _{1\leq i\leq n-1,\ 1\leq j\leq n}$, $B=\left(  X_{i}\right)
_{1\leq i\leq n,\ 1\leq j\leq1}$, $P=\left\{  1,2,\ldots,n\right\}
\setminus\left\{  i\right\}  $ and $Q=\left\{  1,2,\ldots,n\right\}
\setminus\left\{  j\right\}  $ and $\mathbb{K}=k\left[  U_{pq},X_{t}\right]
\diagup\left(  f_{1},f_{2},\ldots,f_{n-1}\right)  $.}:

\begin{exercise}
\label{exe.det.AB=0-proport}Let $n\in\mathbb{N}$. Let $k\in\left\{
0,1,\ldots,n\right\}  $. Let $A\in\mathbb{K}^{k\times n}$ and $B\in
\mathbb{K}^{n\times\left(  n-k\right)  }$ be two matrices satisfying
$AB=0_{k\times\left(  n-k\right)  }$. (Recall that $0_{k\times\left(
n-k\right)  }$ means the $k\times\left(  n-k\right)  $ zero matrix.) Let $P$
and $Q$ be two subsets of $\left\{  1,2,\ldots,n\right\}  $ such that
$\left\vert P\right\vert =\left\vert Q\right\vert =k$.

We shall use the following notations:

\begin{itemize}
\item If $I$ is a finite set of integers, then $w\left(  I\right)  $ shall
denote the list of all elements of $I$ in increasing order (with no
repetitions). (See Definition \ref{def.ind.inclist} for the formal definition
of this list.) (For example, $w\left(  \left\{  3,4,8\right\}  \right)
=\left(  3,4,8\right)  $.)

\item For any subset $I$ of $\left\{  1,2,\ldots,n\right\}  $, we let
$\widetilde{I}$ denote the complement $\left\{  1,2,\ldots,n\right\}
\setminus I$ of $I$. (For instance, if $n=4$ and $I=\left\{  1,4\right\}  $,
then $\widetilde{I}=\left\{  2,3\right\}  $.)

\item If $i_{1},i_{2},\ldots,i_{u}$ are some elements of $\left\{
1,2,\ldots,n\right\}  $, then $\operatorname*{cols}\nolimits_{\left(
i_{1},i_{2},\ldots,i_{u}\right)  }A$ shall denote the matrix
$\operatorname*{cols}\nolimits_{i_{1},i_{2},\ldots,i_{u}}A$.

\item If $i_{1},i_{2},\ldots,i_{u}$ are some elements of $\left\{
1,2,\ldots,n\right\}  $, then $\operatorname*{rows}\nolimits_{\left(
i_{1},i_{2},\ldots,i_{u}\right)  }B$ shall denote the matrix
$\operatorname*{rows}\nolimits_{i_{1},i_{2},\ldots,i_{u}}B$.
\end{itemize}

Prove that%
\begin{align*}
&  \left(  -1\right)  ^{\sum\widetilde{Q}}\det\left(  \operatorname*{cols}%
\nolimits_{w\left(  P\right)  }A\right)  \cdot\det\left(  \operatorname*{rows}%
\nolimits_{w\left(  \widetilde{Q}\right)  }B\right) \\
&  =\left(  -1\right)  ^{\sum\widetilde{P}}\det\left(  \operatorname*{cols}%
\nolimits_{w\left(  Q\right)  }A\right)  \cdot\det\left(  \operatorname*{rows}%
\nolimits_{w\left(  \widetilde{P}\right)  }B\right)  .
\end{align*}

\end{exercise}

\begin{example}
Let us illustrate Exercise \ref{exe.det.AB=0-proport} on an example. Set $n=5$
and $k=2$, and let%
\[
A=\left(
\begin{array}
[c]{ccccc}%
a & b & c & d & e\\
a^{\prime} & b^{\prime} & c^{\prime} & d^{\prime} & e^{\prime}%
\end{array}
\right)  \ \ \ \ \ \ \ \ \ \ \text{and}\ \ \ \ \ \ \ \ \ \ B=\left(
\begin{array}
[c]{ccc}%
x & x^{\prime} & x^{\prime\prime}\\
y & y^{\prime} & y^{\prime\prime}\\
z & z^{\prime} & z^{\prime\prime}\\
u & u^{\prime} & u^{\prime\prime}\\
v & v^{\prime} & v^{\prime\prime}%
\end{array}
\right)
\]
be two matrices satisfying $AB=0_{2\times3}$. Let $P=\left\{  1,3\right\}  $
and $Q=\left\{  3,5\right\}  $. Then, Exercise \ref{exe.det.AB=0-proport}
claims that
\begin{align*}
&  \left(  -1\right)  ^{1+2+4}\det\left(  \operatorname*{cols}\nolimits_{1,3}%
A\right)  \cdot\det\left(  \operatorname*{rows}\nolimits_{1,2,4}B\right) \\
&  =\left(  -1\right)  ^{2+4+5}\det\left(  \operatorname*{cols}\nolimits_{3,5}%
A\right)  \cdot\det\left(  \operatorname*{rows}\nolimits_{2,4,5}B\right)  .
\end{align*}
In other words, it claims that%
\begin{align*}
&  \left(  -1\right)  ^{1+2+4}\det\left(
\begin{array}
[c]{cc}%
a & c\\
a^{\prime} & c^{\prime}%
\end{array}
\right)  \cdot\det\left(
\begin{array}
[c]{ccc}%
x & x^{\prime} & x^{\prime\prime}\\
y & y^{\prime} & y^{\prime\prime}\\
u & u^{\prime} & u^{\prime\prime}%
\end{array}
\right) \\
&  =\left(  -1\right)  ^{2+4+5}\det\left(
\begin{array}
[c]{cc}%
c & e\\
c^{\prime} & e^{\prime}%
\end{array}
\right)  \cdot\det\left(
\begin{array}
[c]{ccc}%
y & y^{\prime} & y^{\prime\prime}\\
u & u^{\prime} & u^{\prime\prime}\\
v & v^{\prime} & v^{\prime\prime}%
\end{array}
\right)  .
\end{align*}
It is not at all straightforward to derive this equality from the condition
$AB=0_{2\times3}$.
\end{example}

The next exercise (taken from \cite[\S 1.2.3, Exercise 44]{Knuth-TAoCP1})
states some deeper properties of the \textit{Cauchy matrix} $\left(  \dfrac
{1}{x_{i}+y_{j}}\right)  _{1\leq i\leq n,\ 1\leq j\leq n}$ whose determinant
we computed in Exercise \ref{exe.cauchy-det}:

\begin{exercise}
\label{exe.cauchy-border}Let $n\in\mathbb{N}$. Let $x_{1},x_{2},\ldots,x_{n}$
be $n$ elements of $\mathbb{K}$. Let $y_{1},y_{2},\ldots,y_{n}$ be $n$
elements of $\mathbb{K}$. Assume that $x_{i}+y_{j}$ is invertible in
$\mathbb{K}$ for every $\left(  i,j\right)  \in\left\{  1,2,\ldots,n\right\}
^{2}$. Let $C$ be the matrix $\left(  \dfrac{1}{x_{i}+y_{j}}\right)  _{1\leq
i\leq n,\ 1\leq j\leq n}$. Prove the following:

\textbf{(a)} The sum of all entries of the adjugate matrix
$\operatorname*{adj}C$ is
\[
\left(  \sum_{k=1}^{n}x_{k}+\sum_{k=1}^{n}y_{k}\right)  \cdot\frac
{\prod_{1\leq i<j\leq n}\left(  \left(  x_{i}-x_{j}\right)  \left(
y_{i}-y_{j}\right)  \right)  }{\prod_{\left(  i,j\right)  \in\left\{
1,2,\ldots,n\right\}  ^{2}}\left(  x_{i}+y_{j}\right)  }.
\]

\textbf{(b)} If the matrix $C$ is invertible, then the sum of all entries of
its inverse $C^{-1}$ is
\[
\sum_{k=1}^{n}x_{k}+\sum_{k=1}^{n}y_{k}.
\]

\textbf{(c)} Let $u\in\mathbb{K}^{n}$ be the column vector whose all $n$
entries equal $1$. Let $D$ be the $\left(  n+1\right)  \times\left(
n+1\right)  $-matrix $\left(
\begin{array}
[c]{cc}%
C & u\\
u^{T} & 0_{1\times1}%
\end{array}
\right)  $ (where we are using the notation from Definition \ref{def.block2x2}%
). Then,%
\[
\det D=-\left(  \sum_{k=1}^{n}x_{k}+\sum_{k=1}^{n}y_{k}\right)  \cdot
\frac{\prod_{1\leq i<j\leq n}\left(  \left(  x_{i}-x_{j}\right)  \left(
y_{i}-y_{j}\right)  \right)  }{\prod_{\left(  i,j\right)  \in\left\{
1,2,\ldots,n\right\}  ^{2}}\left(  x_{i}+y_{j}\right)  }.
\]

[\textbf{Hint:} To prove part \textbf{(a)}, first show that\textit{ }%
\[
\sum_{i=1}^{n}\ \ \sum_{j=1}^{n}\left(  -1\right)  ^{i+j}a_{i,j}\left(
x_{i}+y_{j}\right)  \det\left(  A_{\sim i,\sim j}\right)  =\left(  \sum
_{k=1}^{n}x_{k}+\sum_{k=1}^{n}y_{k}\right)  \cdot\det A
\]
for any $n\times n$-matrix $A=\left(  a_{i,j}\right)  _{1\leq i\leq n,\ 1\leq
j\leq n}$.]
\end{exercise}

\begin{example}
If $n=3$, then the matrix $D$ in Exercise \ref{exe.cauchy-border} \textbf{(c)}
is%
\[
\left(
\begin{array}
[c]{cccc}%
\dfrac{1}{x_{1}+y_{1}} & \dfrac{1}{x_{1}+y_{2}} & \dfrac{1}{x_{1}+y_{3}} & 1\\
\dfrac{1}{x_{2}+y_{1}} & \dfrac{1}{x_{2}+y_{2}} & \dfrac{1}{x_{2}+y_{3}} & 1\\
\dfrac{1}{x_{3}+y_{1}} & \dfrac{1}{x_{3}+y_{2}} & \dfrac{1}{x_{3}+y_{3}} & 1\\
1 & 1 & 1 & 0
\end{array}
\right)  .
\]

\end{example}

\section{Solutions}

This section contains solutions (or, sometimes, solution sketches) to some of
the exercises in the text, as well as occasional remarks. I do not recommend
reading them before trying to solve the problem on your own.

\subsection{Solution to Exercise \ref{exe.jectivity.pigeons}}

\begin{vershort}
\begin{proof}
[Proof of Lemma \ref{lem.jectivity.pigeon0}.]\textbf{(a)} Let $S$ be a subset
of $U$. We must prove that $\left\vert f\left(  S\right)  \right\vert
\leq\left\vert S\right\vert $.

Let $\left(  s_{1},s_{2},\ldots,s_{k}\right)  $ be a list of all elements of
$S$ (with no repetitions).\footnote{Such a list exists, since the set $S$ is
finite.} Thus, $\left\{  s_{1},s_{2},\ldots,s_{k}\right\}  =S$ and
$k=\left\vert S\right\vert $. Now,%
\[
f\left(  \underbrace{S}_{=\left\{  s_{1},s_{2},\ldots,s_{k}\right\}  }\right)
=f\left(  \left\{  s_{1},s_{2},\ldots,s_{k}\right\}  \right)  =\left\{
f\left(  s_{1}\right)  ,f\left(  s_{2}\right)  ,\ldots,f\left(  s_{k}\right)
\right\}  .
\]
Hence,%
\begin{align*}
\left\vert f\left(  S\right)  \right\vert  &  =\left\vert \left\{  f\left(
s_{1}\right)  ,f\left(  s_{2}\right)  ,\ldots,f\left(  s_{k}\right)  \right\}
\right\vert \\
&  \leq k\ \ \ \ \ \ \ \ \ \ \left(
\begin{array}
[c]{c}%
\text{since there are at most }k\text{ distinct elements}\\
\text{among }f\left(  s_{1}\right)  ,f\left(  s_{2}\right)  ,\ldots,f\left(
s_{k}\right)
\end{array}
\right) \\
&  =\left\vert S\right\vert .
\end{align*}
This proves Lemma \ref{lem.jectivity.pigeon0} \textbf{(a)}.

\textbf{(b)} Assume that $\left\vert f\left(  U\right)  \right\vert
\geq\left\vert U\right\vert $.

Let $u$ and $v$ be two elements of $U$ satisfying $f\left(  u\right)
=f\left(  v\right)  $. We shall prove that $u=v$.

Indeed, assume the contrary. Thus, $u\neq v$. Hence, $u\in U\setminus\left\{
v\right\}  $. But $v\in U$ and therefore $\left\vert U\setminus\left\{
v\right\}  \right\vert =\left\vert U\right\vert -1$.

Lemma \ref{lem.jectivity.pigeon0} \textbf{(a)} (applied to $S=U\setminus
\left\{  v\right\}  $) shows that $\left\vert f\left(  U\setminus\left\{
v\right\}  \right)  \right\vert \leq\left\vert U\setminus\left\{  v\right\}
\right\vert =\left\vert U\right\vert -1$.

Next, I claim that $q\in f\left(  U\setminus\left\{  v\right\}  \right)  $ for
each $q\in f\left(  U\right)  $.

Indeed, let $q\in f\left(  U\right)  $ be arbitrary. We want to show that
$q\in f\left(  U\setminus\left\{  v\right\}  \right)  $.

We know that $q\in f\left(  U\right)  $. Hence, there exists some $p\in U$
satisfying $q=f\left(  p\right)  $. Consider this $p$. If $p=v$, then%
\[
q=f\left(  \underbrace{p}_{=v}\right)  =f\left(  v\right)  =f\left(
\underbrace{u}_{\in U\setminus\left\{  v\right\}  }\right)  \in f\left(
U\setminus\left\{  v\right\}  \right)  .
\]
Hence, if $p=v$, then $q\in f\left(  U\setminus\left\{  v\right\}  \right)  $
is proven. Thus, for the rest of the proof of $q\in f\left(  U\setminus
\left\{  v\right\}  \right)  $, we WLOG assume that $p\neq v$.

Hence, $p\in U\setminus\left\{  v\right\}  $. Now, $q=f\left(  \underbrace{p}%
_{\in U\setminus\left\{  v\right\}  }\right)  \in f\left(  U\setminus\left\{
v\right\}  \right)  $. This completes the proof of $q\in f\left(
U\setminus\left\{  v\right\}  \right)  $.

Now, forget that we fixed $q$. We thus have shown that $q\in f\left(
U\setminus\left\{  v\right\}  \right)  $ for each $q\in f\left(  U\right)  $.
In other words, $f\left(  U\right)  \subseteq f\left(  U\setminus\left\{
v\right\}  \right)  $. Hence, $\left\vert f\left(  U\right)  \right\vert
\leq\left\vert f\left(  U\setminus\left\{  v\right\}  \right)  \right\vert
\leq\left\vert U\right\vert -1<\left\vert U\right\vert $. Thus, $\left\vert
U\right\vert >\left\vert f\left(  U\right)  \right\vert \geq\left\vert
U\right\vert $. This is absurd. This contradiction shows that our assumption
was false. Thus, $u=v$ is proven.

Now, forget that we fixed $u$ and $v$. We thus have shown that if $u$ and $v$
are two elements of $U$ satisfying $f\left(  u\right)  =f\left(  v\right)  $,
then $u=v$. In other words, the map $f$ is injective. This proves Lemma
\ref{lem.jectivity.pigeon0} \textbf{(b)}.

\textbf{(c)} Assume that $f$ is injective. Let $S$ be a subset of $U$. We must
prove that $\left\vert f\left(  S\right)  \right\vert =\left\vert S\right\vert
$.

Let $\left(  s_{1},s_{2},\ldots,s_{k}\right)  $ be a list of all elements of
$S$ (with no repetitions).\footnote{Such a list exists, since the set $S$ is
finite.} Thus, $\left\{  s_{1},s_{2},\ldots,s_{k}\right\}  =S$ and
$k=\left\vert S\right\vert $. Furthermore, the elements $s_{1},s_{2}%
,\ldots,s_{k}$ are pairwise distinct (since $\left(  s_{1},s_{2},\ldots
,s_{k}\right)  $ is a list with no repetitions). In other words,%
\[
s_{i}\neq s_{j}\ \ \ \ \ \ \ \ \ \ \text{for any two distinct elements
}i\text{ and }j\text{ of }\left\{  1,2,\ldots,k\right\}  .
\]
Therefore,%
\[
f\left(  s_{i}\right)  \neq f\left(  s_{j}\right)
\ \ \ \ \ \ \ \ \ \ \text{for any two distinct elements }i\text{ and }j\text{
of }\left\{  1,2,\ldots,k\right\}
\]
(since the map $f$ is injective, and thus $f\left(  s_{i}\right)  \neq
f\left(  s_{j}\right)  $ follows from $s_{i}\neq s_{j}$). In other words, the
$k$ elements $f\left(  s_{1}\right)  ,f\left(  s_{2}\right)  ,\ldots,f\left(
s_{k}\right)  $ are pairwise distinct. Hence, $\left\vert \left\{  f\left(
s_{1}\right)  ,f\left(  s_{2}\right)  ,\ldots,f\left(  s_{k}\right)  \right\}
\right\vert =k$.

Now,%
\[
f\left(  \underbrace{S}_{=\left\{  s_{1},s_{2},\ldots,s_{k}\right\}  }\right)
=f\left(  \left\{  s_{1},s_{2},\ldots,s_{k}\right\}  \right)  =\left\{
f\left(  s_{1}\right)  ,f\left(  s_{2}\right)  ,\ldots,f\left(  s_{k}\right)
\right\}  .
\]
Hence,%
\[
\left\vert \underbrace{f\left(  S\right)  }_{=\left\{  f\left(  s_{1}\right)
,f\left(  s_{2}\right)  ,\ldots,f\left(  s_{k}\right)  \right\}  }\right\vert
=\left\vert \left\{  f\left(  s_{1}\right)  ,f\left(  s_{2}\right)
,\ldots,f\left(  s_{k}\right)  \right\}  \right\vert =k=\left\vert
S\right\vert .
\]
This proves Lemma \ref{lem.jectivity.pigeon0} \textbf{(c)}.
\end{proof}
\end{vershort}

\begin{verlong}
\begin{proof}
[Proof of Lemma \ref{lem.jectivity.pigeon0}.]\textbf{(a)} Let $S$ be a subset
of $U$. We must prove that $\left\vert f\left(  S\right)  \right\vert
\leq\left\vert S\right\vert $.

The set $S$ is finite (since it is a subset of the finite set $U$). Let
$\left(  s_{1},s_{2},\ldots,s_{k}\right)  $ be a list of all elements of $S$
(with no repetitions).\footnote{Such a list exists, since the set $S$ is
finite.} Thus, $\left\{  s_{1},s_{2},\ldots,s_{k}\right\}  =S$ and
$k=\left\vert S\right\vert $. Now,%
\begin{equation}
f\left(  \underbrace{S}_{=\left\{  s_{1},s_{2},\ldots,s_{k}\right\}  }\right)
=f\left(  \left\{  s_{1},s_{2},\ldots,s_{k}\right\}  \right)  =\left\{
f\left(  s_{1}\right)  ,f\left(  s_{2}\right)  ,\ldots,f\left(  s_{k}\right)
\right\}  . \label{pf.lem.jectivity.pigeon0.a.1}%
\end{equation}

But the elements of the set $\left\{  f\left(  s_{1}\right)  ,f\left(
s_{2}\right)  ,\ldots,f\left(  s_{k}\right)  \right\}  $ are the $k$ elements
\newline$f\left(  s_{1}\right)  ,f\left(  s_{2}\right)  ,\ldots,f\left(
s_{k}\right)  $, which may and may not be distinct; in either case, there are
at most $k$ distinct elements among them. Hence, the set $\left\{  f\left(
s_{1}\right)  ,f\left(  s_{2}\right)  ,\ldots,f\left(  s_{k}\right)  \right\}
$ has at most $k$ elements. In other words, $\left\vert \left\{  f\left(
s_{1}\right)  ,f\left(  s_{2}\right)  ,\ldots,f\left(  s_{k}\right)  \right\}
\right\vert \leq k$. In light of (\ref{pf.lem.jectivity.pigeon0.a.1}), this
rewrites as $\left\vert f\left(  S\right)  \right\vert \leq k$. Thus,
$\left\vert f\left(  S\right)  \right\vert \leq k=\left\vert S\right\vert $.
This proves Lemma \ref{lem.jectivity.pigeon0} \textbf{(a)}.

\textbf{(b)} Assume that $\left\vert f\left(  U\right)  \right\vert
\geq\left\vert U\right\vert $.

Let $u$ and $v$ be two elements of $U$ satisfying $f\left(  u\right)
=f\left(  v\right)  $. We shall prove that $u=v$.

Indeed, assume the contrary. Thus, $u\neq v$. Hence, $u\in U\setminus\left\{
v\right\}  $ (since $u\in U$ and $u\neq v$). But $v\in U$ and therefore
$\left\vert U\setminus\left\{  v\right\}  \right\vert =\left\vert U\right\vert
-1$.

Lemma \ref{lem.jectivity.pigeon0} \textbf{(a)} (applied to $S=U\setminus
\left\{  v\right\}  $) shows that $\left\vert f\left(  U\setminus\left\{
v\right\}  \right)  \right\vert \leq\left\vert U\setminus\left\{  v\right\}
\right\vert =\left\vert U\right\vert -1$.

Next, I claim that $q\in f\left(  U\setminus\left\{  v\right\}  \right)  $ for
each $q\in f\left(  U\right)  $.

Indeed, let $q\in f\left(  U\right)  $ be arbitrary. We want to show that
$q\in f\left(  U\setminus\left\{  v\right\}  \right)  $.

We know that $q\in f\left(  U\right)  $. Hence, there exists some $p\in U$
satisfying $q=f\left(  p\right)  $. Consider this $p$. If $p=v$, then%
\[
q=f\left(  \underbrace{p}_{=v}\right)  =f\left(  v\right)  =f\left(
\underbrace{u}_{\in U\setminus\left\{  v\right\}  }\right)  \in f\left(
U\setminus\left\{  v\right\}  \right)  .
\]
Hence, if $p=v$, then $q\in f\left(  U\setminus\left\{  v\right\}  \right)  $
is proven. Thus, for the rest of the proof of $q\in f\left(  U\setminus
\left\{  v\right\}  \right)  $, we can WLOG assume that we don't have $p=v$.
Assume this.

We have $p\neq v$ (since we don't have $p=v$). Combining this with $p\in U$,
we obtain $p\in U\setminus\left\{  v\right\}  $. Now, $q=f\left(
\underbrace{p}_{\in U\setminus\left\{  v\right\}  }\right)  \in f\left(
U\setminus\left\{  v\right\}  \right)  $. This completes the proof of $q\in
f\left(  U\setminus\left\{  v\right\}  \right)  $.

Now, forget that we fixed $q$. We thus have shown that $q\in f\left(
U\setminus\left\{  v\right\}  \right)  $ for each $q\in f\left(  U\right)  $.
In other words, $f\left(  U\right)  \subseteq f\left(  U\setminus\left\{
v\right\}  \right)  $. Hence, $\left\vert f\left(  U\right)  \right\vert
\leq\left\vert f\left(  U\setminus\left\{  v\right\}  \right)  \right\vert
\leq\left\vert U\right\vert -1<\left\vert U\right\vert $. Thus, $\left\vert
U\right\vert >\left\vert f\left(  U\right)  \right\vert \geq\left\vert
U\right\vert $. This is absurd. Hence, we have found a contradiction. This
contradiction shows that our assumption was false. Thus, $u=v$ is proven.

Now, forget that we fixed $u$ and $v$. We thus have shown that if $u$ and $v$
are two elements of $U$ satisfying $f\left(  u\right)  =f\left(  v\right)  $,
then $u=v$. In other words, the map $f$ is injective. This proves Lemma
\ref{lem.jectivity.pigeon0} \textbf{(b)}.

\textbf{(c)} Assume that $f$ is injective. Let $S$ be a subset of $U$. We must
prove that $\left\vert f\left(  S\right)  \right\vert =\left\vert S\right\vert
$.

The set $S$ is finite (since it is a subset of the finite set $U$). Let
$\left(  s_{1},s_{2},\ldots,s_{k}\right)  $ be a list of all elements of $S$
(with no repetitions).\footnote{Such a list exists, since the set $S$ is
finite.} Thus, $\left\{  s_{1},s_{2},\ldots,s_{k}\right\}  =S$ and
$k=\left\vert S\right\vert $. Furthermore, the elements $s_{1},s_{2}%
,\ldots,s_{k}$ are pairwise distinct (since $\left(  s_{1},s_{2},\ldots
,s_{k}\right)  $ is a list with no repetitions). In other words,%
\begin{equation}
s_{i}\neq s_{j}\ \ \ \ \ \ \ \ \ \ \text{for any two distinct elements
}i\text{ and }j\text{ of }\left\{  1,2,\ldots,k\right\}  .
\label{pf.lem.jectivity.pigeon0.c.1}%
\end{equation}

Now, we have%
\[
f\left(  s_{i}\right)  \neq f\left(  s_{j}\right)
\ \ \ \ \ \ \ \ \ \ \text{for any two distinct elements }i\text{ and }j\text{
of }\left\{  1,2,\ldots,k\right\}
\]
\footnote{\textit{Proof.} Let $i$ and $j$ be two distinct elements of
$\left\{  1,2,\ldots,k\right\}  $. Then, $s_{i}\neq s_{j}$ (by
(\ref{pf.lem.jectivity.pigeon0.c.1})). Assume (for the sake of contradiction)
that $f\left(  s_{i}\right)  =f\left(  s_{j}\right)  $.
\par
Recall that the map $f$ is injective. In other words, any two elements $u$ and
$v$ of $U$ satisfying $f\left(  u\right)  =f\left(  v\right)  $ must satisfy
$u=v$. Applying this to $u=s_{i}$ and $v=s_{j}$, we obtain $s_{i}=s_{j}$
(since $f\left(  s_{i}\right)  =f\left(  s_{j}\right)  $). This contradicts
$s_{i}\neq s_{j}$.
\par
This contradiction shows that our assumption (that $f\left(  s_{i}\right)
=f\left(  s_{j}\right)  $) was wrong. Hence, we cannot have $f\left(
s_{i}\right)  =f\left(  s_{j}\right)  $. In other words, we have $f\left(
s_{i}\right)  \neq f\left(  s_{j}\right)  $. Qed.}. In other words, the $k$
elements $f\left(  s_{1}\right)  ,f\left(  s_{2}\right)  ,\ldots,f\left(
s_{k}\right)  $ are pairwise distinct. Hence, $\left\vert \left\{  f\left(
s_{1}\right)  ,f\left(  s_{2}\right)  ,\ldots,f\left(  s_{k}\right)  \right\}
\right\vert =k$.

Now,%
\[
f\left(  \underbrace{S}_{=\left\{  s_{1},s_{2},\ldots,s_{k}\right\}  }\right)
=f\left(  \left\{  s_{1},s_{2},\ldots,s_{k}\right\}  \right)  =\left\{
f\left(  s_{1}\right)  ,f\left(  s_{2}\right)  ,\ldots,f\left(  s_{k}\right)
\right\}  .
\]
Hence,%
\[
\left\vert \underbrace{f\left(  S\right)  }_{=\left\{  f\left(  s_{1}\right)
,f\left(  s_{2}\right)  ,\ldots,f\left(  s_{k}\right)  \right\}  }\right\vert
=\left\vert \left\{  f\left(  s_{1}\right)  ,f\left(  s_{2}\right)
,\ldots,f\left(  s_{k}\right)  \right\}  \right\vert =k=\left\vert
S\right\vert .
\]
This proves Lemma \ref{lem.jectivity.pigeon0} \textbf{(c)}.
\end{proof}
\end{verlong}

\begin{vershort}
\begin{proof}
[Proof of Lemma \ref{lem.jectivity.pigeon-surj}.]Assume that $f$ is
surjective. Thus, $f\left(  U\right)  =V$. Hence, $\left\vert f\left(
U\right)  \right\vert =\left\vert V\right\vert \geq\left\vert U\right\vert $
(since $\left\vert U\right\vert \leq\left\vert V\right\vert $). Hence, Lemma
\ref{lem.jectivity.pigeon0} \textbf{(b)} shows that the map $f$ is injective.
Since $f$ is both surjective and injective, we see that $f$ is bijective.

Now, forget that we have assumed that $f$ is surjective. We thus have shown
that if $f$ is surjective, then $f$ is bijective. Of course, the converse also
holds: If $f$ is bijective, then $f$ is surjective. Hence, $f$ is surjective
if and only if $f$ is bijective. This proves Lemma
\ref{lem.jectivity.pigeon-surj}.
\end{proof}
\end{vershort}

\begin{verlong}
\begin{proof}
[Proof of Lemma \ref{lem.jectivity.pigeon-surj}.]Assume that $f$ is
surjective. Thus, $f\left(  U\right)  =V$. Hence, $\left\vert f\left(
U\right)  \right\vert =\left\vert V\right\vert \geq\left\vert U\right\vert $
(since $\left\vert U\right\vert \leq\left\vert V\right\vert $). Hence, Lemma
\ref{lem.jectivity.pigeon0} \textbf{(b)} shows that the map $f$ is injective.
Since $f$ is both surjective and injective, we see that $f$ is bijective.

Now, forget that we have assumed that $f$ is surjective. We thus have shown
that if $f$ is surjective, then $f$ is bijective. In other words, we have
proven the following implication:%
\begin{equation}
\left(  f\text{ is surjective}\right)  \ \Longrightarrow\ \left(  f\text{ is
bijective}\right)  . \label{sol.noncomm.polarization.c.fact1.1}%
\end{equation}
On the other hand, each bijective map is surjective. Hence, we have the
following implication:%
\[
\left(  f\text{ is bijective}\right)  \ \Longrightarrow\ \left(  f\text{ is
surjective}\right)  .
\]
Combining this implication with (\ref{sol.noncomm.polarization.c.fact1.1}), we
obtain the logical equivalence%
\[
\left(  f\text{ is surjective}\right)  \ \Longleftrightarrow\ \left(  f\text{
is bijective}\right)  .
\]
This proves Lemma \ref{lem.jectivity.pigeon-surj}.
\end{proof}
\end{verlong}

\begin{vershort}
\begin{proof}
[Proof of Lemma \ref{lem.jectivity.pigeon-inj}.]Assume that $f$ is injective.
Thus, Lemma \ref{lem.jectivity.pigeon0} \textbf{(c)} (applied to $S=U$) yields
$\left\vert f\left(  U\right)  \right\vert =\left\vert U\right\vert $ (since
$U$ is a subset of $U$). Thus, $\left\vert f\left(  U\right)  \right\vert
=\left\vert U\right\vert \geq\left\vert V\right\vert $.

But the following fact is well-known: If $P$ is a finite set, and if $Q$ is a
subset of $P$ such that $\left\vert Q\right\vert \geq\left\vert P\right\vert
$, then $Q=P$. Applying this to $P=V$ and $Q=f\left(  U\right)  $, we conclude
that $f\left(  U\right)  =V$ (since $f\left(  U\right)  $ is a subset of $V$
such that $\left\vert f\left(  U\right)  \right\vert \geq\left\vert
V\right\vert $). In other words, the map $f$ is surjective. Since $f$ is both
surjective and injective, we see that $f$ is bijective.

Now, forget that we have assumed that $f$ is injective. We thus have shown
that if $f$ is injective, then $f$ is bijective. Of course, the converse also
holds: If $f$ is bijective, then $f$ is injective. Hence, $f$ is injective if
and only if $f$ is bijective. This proves Lemma \ref{lem.jectivity.pigeon-inj}.
\end{proof}
\end{vershort}

\begin{verlong}
\begin{proof}
[Proof of Lemma \ref{lem.jectivity.pigeon-inj}.]Assume that $f$ is injective.
Thus, Lemma \ref{lem.jectivity.pigeon0} \textbf{(c)} (applied to $S=U$) yields
$\left\vert f\left(  U\right)  \right\vert =\left\vert U\right\vert $ (since
$U$ is a subset of $U$). Thus, $\left\vert f\left(  U\right)  \right\vert
=\left\vert U\right\vert \geq\left\vert V\right\vert $.

But the following fact is well-known: If $P$ is a finite set, and if $Q$ is a
subset of $P$ such that $\left\vert Q\right\vert \geq\left\vert P\right\vert
$, then $Q=P$. Applying this to $P=V$ and $Q=f\left(  U\right)  $, we conclude
that $f\left(  U\right)  =V$ (since $f\left(  U\right)  $ is a subset of $V$
such that $\left\vert f\left(  U\right)  \right\vert \geq\left\vert
V\right\vert $). In other words, the map $f$ is surjective. Since $f$ is both
surjective and injective, we see that $f$ is bijective.

Now, forget that we have assumed that $f$ is injective. We thus have shown
that if $f$ is injective, then $f$ is bijective. In other words, we have
proven the following implication:%
\begin{equation}
\left(  f\text{ is injective}\right)  \ \Longrightarrow\ \left(  f\text{ is
bijective}\right)  . \label{pf.lem.jectivity.pigeon-inj.1}%
\end{equation}
On the other hand, each bijective map is injective. Hence, we have the
following implication:%
\[
\left(  f\text{ is bijective}\right)  \ \Longrightarrow\ \left(  f\text{ is
injective}\right)  .
\]
Combining this implication with (\ref{pf.lem.jectivity.pigeon-inj.1}), we
obtain the logical equivalence%
\[
\left(  f\text{ is injective}\right)  \ \Longleftrightarrow\ \left(  f\text{
is bijective}\right)  .
\]
This proves Lemma \ref{lem.jectivity.pigeon-inj}.
\end{proof}
\end{verlong}

\begin{proof}
[Solution to Exercise \ref{exe.jectivity.pigeons}.]We have proven Lemma
\ref{lem.jectivity.pigeon0}, Lemma \ref{lem.jectivity.pigeon-surj} and Lemma
\ref{lem.jectivity.pigeon-inj}. Thus, Exercise \ref{exe.jectivity.pigeons} is solved.
\end{proof}

\subsection{Solution to Exercise \ref{exe.ind.LP2q}}

\begin{proof}
[Solution to Exercise \ref{exe.ind.LP2q}.]First, we notice that the recursive
definition of the sequence $\left(  b_{0},b_{1},b_{2},\ldots\right)  $ yields%
\begin{align*}
b_{2}  &  =\dfrac{b_{2-1}^{2}+q}{b_{2-2}}=\dfrac{b_{1}^{2}+q}{b_{0}}%
=\dfrac{1^{2}+q}{1}\ \ \ \ \ \ \ \ \ \ \left(  \text{since }b_{0}=1\text{ and
}b_{1}=1\right) \\
&  =q+1.
\end{align*}
Comparing this with%
\[
\left(  q+2\right)  \underbrace{b_{2-1}}_{=b_{1}=1}-\underbrace{b_{2-2}%
}_{=b_{0}=1}=\left(  q+2\right)  \cdot1-1=q+1,
\]
we obtain $b_{2}=\left(  q+2\right)  b_{2-1}-b_{2-2}$.

\textbf{(a)} We shall prove Exercise \ref{exe.ind.LP2q} \textbf{(a)} by
induction on $n$ starting at $2$:

\textit{Induction base:} We have already shown that $b_{2}=\left(  q+2\right)
b_{2-1}-b_{2-2}$. In other words, Exercise \ref{exe.ind.LP2q} \textbf{(a)}
holds for $n=2$. This completes the induction base.

\textit{Induction step:} Let $m\in\mathbb{Z}_{\geq2}$. Assume that Exercise
\ref{exe.ind.LP2q} \textbf{(a)} holds for $n=m$. We must prove that Exercise
\ref{exe.ind.LP2q} \textbf{(a)} holds for $n=m+1$.

We have assumed that Exercise \ref{exe.ind.LP2q} \textbf{(a)} holds for $n=m$.
In other words, we have $b_{m}=\left(  q+2\right)  b_{m-1}-b_{m-2}$. Thus,%
\begin{equation}
b_{m}-\left(  q+2\right)  b_{m-1}=-b_{m-2}. \label{sol.exe.ind.LP2q.a.0}%
\end{equation}

We have $m\in\mathbb{Z}_{\geq2}$. Thus, $m$ is an integer that is $\geq2$.
Hence, the recursive definition of the sequence $\left(  b_{0},b_{1}%
,b_{2},\ldots\right)  $ yields%
\[
b_{m}=\dfrac{b_{m-1}^{2}+q}{b_{m-2}}.
\]
Multiplying this equality by $b_{m-2}$, we obtain
\begin{equation}
b_{m}b_{m-2}=b_{m-1}^{2}+q. \label{sol.exe.ind.LP2q.a.1}%
\end{equation}

Also, $m+1\geq m\geq2$ (since $m$ is $\geq2$). Hence, the recursive definition
of the sequence $\left(  b_{0},b_{1},b_{2},\ldots\right)  $ yields%
\[
b_{m+1}=\dfrac{b_{\left(  m+1\right)  -1}^{2}+q}{b_{\left(  m+1\right)  -2}%
}=\dfrac{b_{m}^{2}+q}{b_{m-1}}.
\]
Multiplying this equality by $b_{m-1}$, we obtain $b_{m-1}b_{m+1}=b_{m}^{2}%
+q$. Now,%
\begin{align*}
b_{m-1}\left(  b_{m+1}-\left(  q+2\right)  b_{m}+b_{m-1}\right)   &
=\underbrace{b_{m-1}b_{m+1}}_{=b_{m}^{2}+q}-\left(  q+2\right)  b_{m-1}%
b_{m}+b_{m-1}^{2}\\
&  =b_{m}^{2}+q-\left(  q+2\right)  b_{m-1}b_{m}+b_{m-1}^{2}\\
&  =b_{m}\underbrace{\left(  b_{m}-\left(  q+2\right)  b_{m-1}\right)
}_{\substack{=-b_{m-2}\\\text{(by (\ref{sol.exe.ind.LP2q.a.0}))}%
}}+\underbrace{b_{m-1}^{2}+q}_{\substack{=b_{m}b_{m-2}\\\text{(by
(\ref{sol.exe.ind.LP2q.a.1}))}}}\\
&  =b_{m}\left(  -b_{m-2}\right)  +b_{m}b_{m-2}=0.
\end{align*}
We can cancel $b_{m-1}$ from this equality (since $b_{m-1}\neq0$ (because
$b_{m-1}$ is a positive rational number)). Thus, we obtain $b_{m+1}-\left(
q+2\right)  b_{m}+b_{m-1}=0$. Hence,%
\[
b_{m+1}=\left(  q+2\right)  b_{m}-b_{m-1}.
\]
Comparing this with $\left(  q+2\right)  \underbrace{b_{\left(  m+1\right)
-1}}_{=b_{m}}+\underbrace{b_{\left(  m+1\right)  -2}}_{=b_{m-1}}=\left(
q+2\right)  b_{m}-b_{m-1}$, we obtain $b_{m+1}=\left(  q+2\right)  b_{\left(
m+1\right)  -1}+b_{\left(  m+1\right)  -2}$. In other words, Exercise
\ref{exe.ind.LP2q} \textbf{(a)} holds for $n=m+1$. This completes the
induction step. Hence, Exercise \ref{exe.ind.LP2q} \textbf{(a)} is proven by induction.

\textbf{(b)} We shall prove Exercise \ref{exe.ind.LP2q} \textbf{(b)} by strong
induction on $n$ starting at $0$:

\textit{Induction step:} Let $m\in\mathbb{N}$.\ \ \ \ \footnote{In order to
match the notations used in Theorem \ref{thm.ind.SIP}, we should be saying
\textquotedblleft Let $m\in\mathbb{Z}_{\geq0}$\textquotedblright\ here, rather
than \textquotedblleft Let $m\in\mathbb{N}$\textquotedblright. But of course,
this amounts to the same thing, since $\mathbb{N}=\mathbb{Z}_{\geq0}$.} Assume
that Exercise \ref{exe.ind.LP2q} \textbf{(b)} holds for every $n\in\mathbb{N}$
satisfying $n<m$. We must now show that Exercise \ref{exe.ind.LP2q}
\textbf{(b)} holds for $n=m$.

We have assumed that Exercise \ref{exe.ind.LP2q} \textbf{(b)} holds for every
$n\in\mathbb{N}$ satisfying $n<m$. In other words, we have
\begin{equation}
b_{n}\in\mathbb{N}\text{ for every }n\in\mathbb{N}\text{ satisfying }n<m.
\label{sol.ind.LP2q.b.IH}%
\end{equation}

We must now show that Exercise \ref{exe.ind.LP2q} \textbf{(b)} holds for
$n=m$. In other words, we must show that $b_{m}\in\mathbb{N}$.

Recall that $\left(  b_{0},b_{1},b_{2},\ldots\right)  $ is a sequence of
positive rational numbers. Thus, $b_{m}$ is a positive rational number.

We are in one of the following three cases:

\textit{Case 1:} We have $m=0$.

\textit{Case 2:} We have $m=1$.

\textit{Case 3:} We have $m>1$.

Let us first consider Case 1. In this case, we have $m=0$. Thus, $b_{m}%
=b_{0}=1\in\mathbb{N}$. Hence, $b_{m}\in\mathbb{N}$ is proven in Case 1.

Similarly, we can prove $b_{m}\in\mathbb{N}$ in Case 2 (using $b_{1}=1$). It
thus remains to prove $b_{m}\in\mathbb{N}$ in Case 3.

So let us consider Case 3. In this case, we have $m>1$. Thus, $m\geq2$ (since
$m$ is an integer), so that $m\in\mathbb{Z}_{\geq2}$. Thus, Exercise
\ref{exe.ind.LP2q} \textbf{(a)} (applied to $n=m$) yields $b_{m}=\left(
q+2\right)  b_{m-1}-b_{m-2}$.

But $m\geq2\geq1$, so that $m-1\in\mathbb{N}$ and $m-1<m$. Hence,
(\ref{sol.ind.LP2q.b.IH}) (applied to $n=m-1$) yields $b_{m-1}\in
\mathbb{N}\subseteq\mathbb{Z}$.

Also, $m\geq2$, so that $m-2\in\mathbb{N}$ and $m-2<m$. Hence,
(\ref{sol.ind.LP2q.b.IH}) (applied to $n=m-2$) yields $b_{m-2}\in
\mathbb{N}\subseteq\mathbb{Z}$.

So we know that $b_{m-1}$ and $b_{m-2}$ are both integers (since $b_{m-2}%
\in\mathbb{Z}$ and $b_{m-1}\in\mathbb{Z}$). Hence, $\left(  q+2\right)
b_{m-1}-b_{m-2}$ is an integer as well (since $q$ is an integer). In other
words, $b_{m}$ is an integer (because $b_{m}=\left(  q+2\right)
b_{m-1}-b_{m-2}$). Since $b_{m}$ is positive, we thus conclude that $b_{m}$ is
a positive integer. Hence, $b_{m}\in\mathbb{N}$. This shows that $b_{m}%
\in\mathbb{N}$ in Case 3.

We now have proven $b_{m}\in\mathbb{N}$ in each of the three Cases 1, 2 and 3.
Since these three Cases cover all possibilities, we thus conclude that
$b_{m}\in\mathbb{N}$ always holds. In other words, Exercise \ref{exe.ind.LP2q}
\textbf{(b)} holds for $n=m$. This completes the induction step. Thus,
Exercise \ref{exe.ind.LP2q} \textbf{(b)} is proven by strong induction.
\end{proof}

\subsection{Solution to Exercise \ref{exe.ind.inverse-comp}}

\begin{proof}
[Proof of Proposition \ref{prop.ind.inverse-comp}.]We shall prove Proposition
\ref{prop.ind.inverse-comp} by induction on $n$:

\textit{Induction base:} We want to prove that Proposition
\ref{prop.ind.inverse-comp} holds when $n=0$. In other words, we want to prove
the following claim:

\begin{statement}
\textit{Claim 1:} Let $X_{1},X_{2},\ldots,X_{0+1}$ be $0+1$ sets. For each
$i\in\left\{  1,2,\ldots,0\right\}  $, let $f_{i}:X_{i}\rightarrow X_{i+1}$ be
an invertible map. Then, the map $f_{0}\circ f_{0-1}\circ\cdots\circ
f_{1}:X_{1}\rightarrow X_{0+1}$ is invertible as well, and its inverse is%
\[
\left(  f_{0}\circ f_{0-1}\circ\cdots\circ f_{1}\right)  ^{-1}=f_{1}^{-1}\circ
f_{2}^{-1}\circ\cdots\circ f_{0}^{-1}.
\]

\end{statement}

But this is straightforward:

[\textit{Proof of Claim 1:} The equality
(\ref{eq.def.ind.gen-ass-maps.comp0.0}) (applied to $n=0$) yields $f_{0}\circ
f_{0-1}\circ\cdots\circ f_{1}=\operatorname*{id}\nolimits_{X_{1}}$. Similarly,%
\begin{equation}
f_{1}^{-1}\circ f_{2}^{-1}\circ\cdots\circ f_{0}^{-1}=\operatorname*{id}%
\nolimits_{X_{1}}. \label{pf.prop.ind.inverse-comp.base.1}%
\end{equation}

Now, $f_{0}\circ f_{0-1}\circ\cdots\circ f_{1}=\operatorname*{id}%
\nolimits_{X_{1}}$ and $X_{0+1}=X_{1}$. Hence, the map $f_{0}\circ
f_{0-1}\circ\cdots\circ f_{1}:X_{1}\rightarrow X_{0+1}$ is the same as the map
$\operatorname*{id}\nolimits_{X_{1}}:X_{1}\rightarrow X_{1}$, and therefore is
invertible (because the map $\operatorname*{id}\nolimits_{X_{1}}%
:X_{1}\rightarrow X_{1}$ clearly is invertible). Moreover, its inverse is%
\[
\left(  \underbrace{f_{0}\circ f_{0-1}\circ\cdots\circ f_{1}}%
_{=\operatorname*{id}\nolimits_{X_{1}}}\right)  ^{-1}=\left(
\operatorname*{id}\nolimits_{X_{1}}\right)  ^{-1}=\operatorname*{id}%
\nolimits_{X_{1}}=f_{1}^{-1}\circ f_{2}^{-1}\circ\cdots\circ f_{0}%
^{-1}\ \ \ \ \ \ \ \ \ \ \left(  \text{by
(\ref{pf.prop.ind.inverse-comp.base.1})}\right)  .
\]
This completes the proof of Claim 1.]

We have now proven Claim 1. In other words, Proposition
\ref{prop.ind.inverse-comp} holds when $n=0$. This completes the induction base.

\textit{Induction step:} Let $m\in\mathbb{N}$. Assume that Proposition
\ref{prop.ind.inverse-comp} holds when $n=m$. We must now prove that
Proposition \ref{prop.ind.inverse-comp} holds when $n=m+1$.

We have assumed that Proposition \ref{prop.ind.inverse-comp} holds when $n=m$.
In other words, the following claim holds:

\begin{statement}
\textit{Claim 2:} Let $X_{1},X_{2},\ldots,X_{m+1}$ be $m+1$ sets. For each
$i\in\left\{  1,2,\ldots,m\right\}  $, let $f_{i}:X_{i}\rightarrow X_{i+1}$ be
an invertible map. Then, the map $f_{m}\circ f_{m-1}\circ\cdots\circ
f_{1}:X_{1}\rightarrow X_{m+1}$ is invertible as well, and its inverse is%
\[
\left(  f_{m}\circ f_{m-1}\circ\cdots\circ f_{1}\right)  ^{-1}=f_{1}^{-1}\circ
f_{2}^{-1}\circ\cdots\circ f_{m}^{-1}.
\]

\end{statement}

We must now prove that Proposition \ref{prop.ind.inverse-comp} holds when
$n=m+1$. In other words, we must prove the following claim:

\begin{statement}
\textit{Claim 3:} Let $X_{1},X_{2},\ldots,X_{\left(  m+1\right)  +1}$ be
$\left(  m+1\right)  +1$ sets. For each $i\in\left\{  1,2,\ldots,m+1\right\}
$, let $f_{i}:X_{i}\rightarrow X_{i+1}$ be an invertible map. Then, the map
$f_{m+1}\circ f_{\left(  m+1\right)  -1}\circ\cdots\circ f_{1}:X_{1}%
\rightarrow X_{\left(  m+1\right)  +1}$ is invertible as well, and its inverse
is%
\[
\left(  f_{m+1}\circ f_{\left(  m+1\right)  -1}\circ\cdots\circ f_{1}\right)
^{-1}=f_{1}^{-1}\circ f_{2}^{-1}\circ\cdots\circ f_{m+1}^{-1}.
\]

\end{statement}

\begin{vershort}
[\textit{Proof of Claim 3:} We have $m\in\mathbb{N}$; thus, $m+1$ is a
positive integer. Hence, $m+1\in\left\{  1,2,\ldots,m+1\right\}  $ and
$m+1\geq1$. Thus, Theorem \ref{thm.ind.gen-ass-maps.1} \textbf{(b)} (applied
to $n=m+1$) yields%
\begin{align}
f_{m+1}\circ f_{\left(  m+1\right)  -1}\circ\cdots\circ f_{1}  &
=f_{m+1}\circ\underbrace{\left(  f_{\left(  m+1\right)  -1}\circ f_{\left(
m+1\right)  -2}\circ\cdots\circ f_{1}\right)  }_{=f_{m}\circ f_{m-1}%
\circ\cdots\circ f_{1}}\nonumber\\
&  =f_{m+1}\circ\left(  f_{m}\circ f_{m-1}\circ\cdots\circ f_{1}\right)  .
\label{pf.prop.ind.inverse-comp.c3.pf.short.0}%
\end{align}
Also, Theorem \ref{thm.ind.gen-ass-maps.1} \textbf{(c)} (applied to $m+1$,
$X_{m+3-i}$ and $f_{m+2-i}^{-1}$ instead of $n$, $X_{i}$ and $f_{i}$) yields%
\begin{equation}
f_{1}^{-1}\circ f_{2}^{-1}\circ\cdots\circ f_{m+1}^{-1}=\left(  f_{1}%
^{-1}\circ f_{2}^{-1}\circ\cdots\circ f_{m}^{-1}\right)  \circ f_{m+1}^{-1}.
\label{pf.prop.ind.inverse-comp.c3.pf.short.4}%
\end{equation}

Claim 2 yields that the map $f_{m}\circ f_{m-1}\circ\cdots\circ f_{1}%
:X_{1}\rightarrow X_{m+1}$ is invertible, and that its inverse is%
\[
\left(  f_{m}\circ f_{m-1}\circ\cdots\circ f_{1}\right)  ^{-1}=f_{1}^{-1}\circ
f_{2}^{-1}\circ\cdots\circ f_{m}^{-1}.
\]
Hence,%
\begin{align}
&  \underbrace{\left(  f_{m}\circ f_{m-1}\circ\cdots\circ f_{1}\right)  ^{-1}%
}_{=f_{1}^{-1}\circ f_{2}^{-1}\circ\cdots\circ f_{m}^{-1}}\circ f_{m+1}%
^{-1}\nonumber\\
&  =\left(  f_{1}^{-1}\circ f_{2}^{-1}\circ\cdots\circ f_{m}^{-1}\right)
\circ f_{m+1}^{-1}\nonumber\\
&  =f_{1}^{-1}\circ f_{2}^{-1}\circ\cdots\circ f_{m+1}^{-1}%
\ \ \ \ \ \ \ \ \ \ \left(  \text{by
(\ref{pf.prop.ind.inverse-comp.c3.pf.short.4})}\right)  .
\label{pf.prop.ind.inverse-comp.c3.pf.short.7}%
\end{align}

We know (by our assumption) that for each $i\in\left\{  1,2,\ldots
,m+1\right\}  $, the map $f_{i}:X_{i}\rightarrow X_{i+1}$ is an invertible
map. Applying this to $i=m+1$, we conclude that $f_{m+1}:X_{m+1}\rightarrow
X_{\left(  m+1\right)  +1}$ is an invertible map (since $m+1\in\left\{
1,2,\ldots,m+1\right\}  $).

Now, we know that $f_{m}\circ f_{m-1}\circ\cdots\circ f_{1}:X_{1}\rightarrow
X_{m+1}$ and $f_{m+1}:X_{m+1}\rightarrow X_{\left(  m+1\right)  +1}$ are two
invertible maps. Hence, Proposition \ref{prop.ind.inverse-fg} (applied to
$X=X_{1}$, $Y=X_{m+1}$, $Z=X_{\left(  m+1\right)  +1}$, $b=f_{m}\circ
f_{m-1}\circ\cdots\circ f_{1}$ and $a=f_{m+1}$) yields that the map
$f_{m+1}\circ\left(  f_{m}\circ f_{m-1}\circ\cdots\circ f_{1}\right)
:X_{1}\rightarrow X_{\left(  m+1\right)  +1}$ is invertible as well, and that
its inverse is%
\[
\left(  f_{m+1}\circ\left(  f_{m}\circ f_{m-1}\circ\cdots\circ f_{1}\right)
\right)  ^{-1}=\left(  f_{m}\circ f_{m-1}\circ\cdots\circ f_{1}\right)
^{-1}\circ f_{m+1}^{-1}.
\]
In view of
\[
f_{m+1}\circ\left(  f_{m}\circ f_{m-1}\circ\cdots\circ f_{1}\right)
=f_{m+1}\circ f_{\left(  m+1\right)  -1}\circ\cdots\circ f_{1}%
\ \ \ \ \ \ \ \ \ \ \left(  \text{by
(\ref{pf.prop.ind.inverse-comp.c3.pf.short.0})}\right)
\]
and%
\[
\left(  f_{m}\circ f_{m-1}\circ\cdots\circ f_{1}\right)  ^{-1}\circ
f_{m+1}^{-1}=f_{1}^{-1}\circ f_{2}^{-1}\circ\cdots\circ f_{m+1}^{-1}%
\ \ \ \ \ \ \ \ \ \ \left(  \text{by
(\ref{pf.prop.ind.inverse-comp.c3.pf.short.7})}\right)  ,
\]
this rewrites as follows: The map $f_{m+1}\circ f_{\left(  m+1\right)
-1}\circ\cdots\circ f_{1}:X_{1}\rightarrow X_{\left(  m+1\right)  +1}$ is
invertible as well, and its inverse is%
\[
\left(  f_{m+1}\circ f_{\left(  m+1\right)  -1}\circ\cdots\circ f_{1}\right)
^{-1}=f_{1}^{-1}\circ f_{2}^{-1}\circ\cdots\circ f_{m+1}^{-1}.
\]
This proves Claim 3.]
\end{vershort}

\begin{verlong}
[\textit{Proof of Claim 3:} We have $m\in\mathbb{N}$; thus, $m+1$ is a
positive integer. Hence, $m+1\in\left\{  1,2,\ldots,m+1\right\}  $ and
$m+1\geq1$. Thus, Theorem \ref{thm.ind.gen-ass-maps.1} \textbf{(b)} (applied
to $n=m+1$) yields%
\begin{align}
f_{m+1}\circ f_{\left(  m+1\right)  -1}\circ\cdots\circ f_{1}  &
=f_{m+1}\circ\underbrace{\left(  f_{\left(  m+1\right)  -1}\circ f_{\left(
m+1\right)  -2}\circ\cdots\circ f_{1}\right)  }_{\substack{=f_{m}\circ
f_{m-1}\circ\cdots\circ f_{1}\\\text{(since }\left(  m+1\right)  -1=m\text{
and }\left(  m+1\right)  -2=m-1\text{)}}}\nonumber\\
&  =f_{m+1}\circ\left(  f_{m}\circ f_{m-1}\circ\cdots\circ f_{1}\right)  .
\label{pf.prop.ind.inverse-comp.c3.pf.0}%
\end{align}

Also, we know that
\begin{equation}
f_{i}:X_{i}\rightarrow X_{i+1}\text{ is an invertible map}
\label{pf.prop.ind.inverse-comp.c3.pf.1}%
\end{equation}
for each $i\in\left\{  1,2,\ldots,m+1\right\}  $. Hence,%
\begin{equation}
f_{m+2-i}^{-1}:X_{m+3-i}\rightarrow X_{m+3-\left(  i+1\right)  }\text{ is a
map} \label{pf.prop.ind.inverse-comp.c3.pf.1rev}%
\end{equation}
for each $i\in\left\{  1,2,\ldots,m+1\right\}  $%
\ \ \ \ \footnote{\textit{Proof of (\ref{pf.prop.ind.inverse-comp.c3.pf.1rev}%
):} Let $i\in\left\{  1,2,\ldots,m+1\right\}  $. Thus, $m+2-i\in\left\{
\left(  m+2\right)  -\left(  m+1\right)  ,\left(  m+2\right)  -m,\ldots
,\left(  m+2\right)  -1\right\}  =\left\{  1,2,\ldots,m+1\right\}  $. Hence,
(\ref{pf.prop.ind.inverse-comp.c3.pf.1}) (applied to $m+2-i$ instead of $i$)
yields that $f_{m+2-i}:X_{m+2-i}\rightarrow X_{m+2-i+1}$ is an invertible map.
Thus, its inverse $f_{m+2-i}^{-1}$ is a map from $X_{m+2-i+1}$ to $X_{m+2-i}$.
In view of $X_{m+2-i+1}=X_{m+3-i}$ (since $m+2-i+1=m+3-i$) and $X_{m+2-i}%
=X_{m+3-\left(  i+1\right)  }$ (since $m+2-i=m+3-\left(  i+1\right)  $), this
rewrites as follows: $f_{m+2-i}^{-1}$ is a map from $X_{m+3-i}$ to
$X_{m+3-\left(  i+1\right)  }$. In other words, $f_{m+2-i}^{-1}:X_{m+3-i}%
\rightarrow X_{m+3-\left(  i+1\right)  }$ is a map. This proves
(\ref{pf.prop.ind.inverse-comp.c3.pf.1rev}).}. Therefore, Theorem
\ref{thm.ind.gen-ass-maps.1} \textbf{(c)} (applied to $m+1$, $X_{m+3-i}$ and
$f_{m+2-i}^{-1}$ instead of $n$, $X_{i}$ and $f_{i}$) yields%
\begin{align*}
&  f_{m+2-\left(  m+1\right)  }^{-1}\circ f_{m+2-\left(  \left(  m+1\right)
-1\right)  }^{-1}\circ\cdots\circ f_{m+2-1}^{-1}\\
&  =\left(  f_{m+2-\left(  m+1\right)  }^{-1}\circ f_{m+2-\left(  \left(
m+1\right)  -1\right)  }^{-1}\circ\cdots\circ f_{m+2-2}^{-1}\right)  \circ
f_{m+2-1}^{-1}.
\end{align*}
In view of $m+2-\left(  m+1\right)  =1$, $m+2-\left(  \left(  m+1\right)
-1\right)  =2$, $m+2-2=m$ and $m+2-1=m+1$, this rewrites as follows:%
\begin{equation}
f_{1}^{-1}\circ f_{2}^{-1}\circ\cdots\circ f_{m+1}^{-1}=\left(  f_{1}%
^{-1}\circ f_{2}^{-1}\circ\cdots\circ f_{m}^{-1}\right)  \circ f_{m+1}^{-1}.
\label{pf.prop.ind.inverse-comp.c3.pf.4}%
\end{equation}

For each $i\in\left\{  1,2,\ldots,m\right\}  $, the map $f_{i}:X_{i}%
\rightarrow X_{i+1}$ is an invertible map (by
(\ref{pf.prop.ind.inverse-comp.c3.pf.1}), because $i\in\left\{  1,2,\ldots
,m\right\}  \subseteq\left\{  1,2,\ldots,m+1\right\}  $). Hence, Claim 2
yields that the map $f_{m}\circ f_{m-1}\circ\cdots\circ f_{1}:X_{1}\rightarrow
X_{m+1}$ is invertible as well, and that its inverse is%
\[
\left(  f_{m}\circ f_{m-1}\circ\cdots\circ f_{1}\right)  ^{-1}=f_{1}^{-1}\circ
f_{2}^{-1}\circ\cdots\circ f_{m}^{-1}.
\]
Hence,%
\begin{align}
&  \underbrace{\left(  f_{m}\circ f_{m-1}\circ\cdots\circ f_{1}\right)  ^{-1}%
}_{=f_{1}^{-1}\circ f_{2}^{-1}\circ\cdots\circ f_{m}^{-1}}\circ f_{m+1}%
^{-1}\nonumber\\
&  =\left(  f_{1}^{-1}\circ f_{2}^{-1}\circ\cdots\circ f_{m}^{-1}\right)
\circ f_{m+1}^{-1}\nonumber\\
&  =f_{1}^{-1}\circ f_{2}^{-1}\circ\cdots\circ f_{m+1}^{-1}%
\ \ \ \ \ \ \ \ \ \ \left(  \text{by (\ref{pf.prop.ind.inverse-comp.c3.pf.4}%
)}\right)  . \label{pf.prop.ind.inverse-comp.c3.pf.7}%
\end{align}

Applying (\ref{pf.prop.ind.inverse-comp.c3.pf.1}) to $i=m+1$, we conclude that
$f_{m+1}:X_{m+1}\rightarrow X_{\left(  m+1\right)  +1}$ is an invertible map
(since $m+1\in\left\{  1,2,\ldots,m+1\right\}  $).

Now, we know that $f_{m}\circ f_{m-1}\circ\cdots\circ f_{1}:X_{1}\rightarrow
X_{m+1}$ and $f_{m+1}:X_{m+1}\rightarrow X_{\left(  m+1\right)  +1}$ are two
invertible maps. Hence, Proposition \ref{prop.ind.inverse-fg} (applied to
$X=X_{1}$, $Y=X_{m+1}$, $Z=X_{\left(  m+1\right)  +1}$, $b=f_{m}\circ
f_{m-1}\circ\cdots\circ f_{1}$ and $a=f_{m+1}$) yields that the map
$f_{m+1}\circ\left(  f_{m}\circ f_{m-1}\circ\cdots\circ f_{1}\right)
:X_{1}\rightarrow X_{\left(  m+1\right)  +1}$ is invertible as well, and that
its inverse is%
\[
\left(  f_{m+1}\circ\left(  f_{m}\circ f_{m-1}\circ\cdots\circ f_{1}\right)
\right)  ^{-1}=\left(  f_{m}\circ f_{m-1}\circ\cdots\circ f_{1}\right)
^{-1}\circ f_{m+1}^{-1}.
\]

In view of
\[
f_{m+1}\circ\left(  f_{m}\circ f_{m-1}\circ\cdots\circ f_{1}\right)
=f_{m+1}\circ f_{\left(  m+1\right)  -1}\circ\cdots\circ f_{1}%
\ \ \ \ \ \ \ \ \ \ \left(  \text{by (\ref{pf.prop.ind.inverse-comp.c3.pf.0}%
)}\right)
\]
and%
\[
\left(  f_{m}\circ f_{m-1}\circ\cdots\circ f_{1}\right)  ^{-1}\circ
f_{m+1}^{-1}=f_{1}^{-1}\circ f_{2}^{-1}\circ\cdots\circ f_{m+1}^{-1}%
\ \ \ \ \ \ \ \ \ \ \left(  \text{by (\ref{pf.prop.ind.inverse-comp.c3.pf.7}%
)}\right)  ,
\]
this rewrites as follows: The map $f_{m+1}\circ f_{\left(  m+1\right)
-1}\circ\cdots\circ f_{1}:X_{1}\rightarrow X_{\left(  m+1\right)  +1}$ is
invertible as well, and its inverse is%
\[
\left(  f_{m+1}\circ f_{\left(  m+1\right)  -1}\circ\cdots\circ f_{1}\right)
^{-1}=f_{1}^{-1}\circ f_{2}^{-1}\circ\cdots\circ f_{m+1}^{-1}.
\]
This proves Claim 3.]
\end{verlong}

We have now proven Claim 3. In other words, Proposition
\ref{prop.ind.inverse-comp} holds when $n=m+1$. This completes the induction
step. Hence, Proposition \ref{prop.ind.inverse-comp} is proven by induction.
\end{proof}

\begin{proof}
[Solution to Exercise \ref{exe.ind.inverse-comp}.]We have proven Proposition
\ref{prop.ind.inverse-comp}. Thus, Exercise \ref{exe.ind.inverse-comp} is solved.
\end{proof}

\subsection{Solution to Exercise \ref{exe.ind.gen-com.fin-sup.proofs}}

Let us first prove a lemma:

\begin{lemma}
\label{lem.sol.ind.gen-com.fin-sup.proofs.1}Let $S$ be any set. Let $\left(
a_{s}\right)  _{s\in S}$ be an $\mathbb{A}$-valued $S$-family. Let $T_{1}$ and
$T_{2}$ be two finite subsets $T$ of $S$ such that
(\ref{eq.def.ind.gen-com.fin-sup.sum.0}) holds. Then,
\[
\sum_{s\in T_{1}}a_{s}=\sum_{s\in T_{2}}a_{s}.
\]

\end{lemma}

(Note that we are not using the notation introduced in Definition
\ref{def.ind.gen-com.fin-sup.sum} yet, because we have not proven that this
notation is well-defined.)

\begin{proof}
[Proof of Lemma \ref{lem.sol.ind.gen-com.fin-sup.proofs.1}.]We know that
$T_{1}$ is a finite subset $T$ of $S$ such that
(\ref{eq.def.ind.gen-com.fin-sup.sum.0}) holds. Hence, $T_{1}$ is a finite
subset of $S$.

We know that $T_{2}$ is a finite subset $T$ of $S$ such that
(\ref{eq.def.ind.gen-com.fin-sup.sum.0}) holds. In other words, $T_{2}$ is a
finite subset of $S$ and has the property that
(\ref{eq.def.ind.gen-com.fin-sup.sum.0}) holds for $T=T_{2}$.

In particular, (\ref{eq.def.ind.gen-com.fin-sup.sum.0}) holds for $T=T_{2}$.
In other words,%
\begin{equation}
\text{every }s\in S\setminus T_{2}\text{ satisfies }a_{s}=0.
\label{pf.lem.sol.ind.gen-com.fin-sup.proofs.1.2}%
\end{equation}

The set $T_{1}\setminus T_{2}$ is a subset of the finite set $T_{1}$, and thus
is finite.

If $s\in T_{1}\setminus T_{2}$, then $s\in\underbrace{T_{1}}_{\subseteq
S}\setminus T_{2}\subseteq S\setminus T_{2}$ and thus
\begin{equation}
a_{s}=0 \label{pf.lem.sol.ind.gen-com.fin-sup.proofs.1.3}%
\end{equation}
(by (\ref{pf.lem.sol.ind.gen-com.fin-sup.proofs.1.2})). Now,%
\[
\sum_{s\in T_{1}\setminus T_{2}}\underbrace{a_{s}}_{\substack{=0\\\text{(by
(\ref{pf.lem.sol.ind.gen-com.fin-sup.proofs.1.3}))}}}=\sum_{s\in
T_{1}\setminus T_{2}}0=0
\]
(by Theorem \ref{thm.ind.gen-com.sum(0)} (applied to $T_{1}\setminus T_{2}$
instead of $S$)).

\begin{vershort}
It is a straightforward exercise in set theory to show that if $A$ and $B$ are
any two sets, then $A\cap B$ and $A\setminus B$ are two subsets of $A$
satisfying
\[
\left(  A\cap B\right)  \cap\left(  A\setminus B\right)  =\varnothing
\ \ \ \ \ \ \ \ \ \ \text{and}\ \ \ \ \ \ \ \ \ \ \left(  A\cap B\right)
\cup\left(  A\setminus B\right)  =A.
\]
Applying this to $A=T_{1}$ and $B=T_{2}$, we conclude that $T_{1}\cap T_{2}$
and $T_{1}\setminus T_{2}$ are two subsets of $T_{1}$ satisfying
\[
\left(  T_{1}\cap T_{2}\right)  \cap\left(  T_{1}\setminus T_{2}\right)
=\varnothing\ \ \ \ \ \ \ \ \ \ \text{and}\ \ \ \ \ \ \ \ \ \ \left(
T_{1}\cap T_{2}\right)  \cup\left(  T_{1}\setminus T_{2}\right)  =T_{1}.
\]

\end{vershort}

\begin{verlong}
Clearly, $T_{1}\cap T_{2}$ and $T_{1}\setminus T_{2}$ are two subsets of
$T_{1}$.

A well-known fact from set theory says that if $A$ is a subset of a set $B$,
then
\[
A\cap\left(  B\setminus A\right)  =\varnothing\ \ \ \ \ \ \ \ \ \ \text{and}%
\ \ \ \ \ \ \ \ \ \ A\cup\left(  B\setminus A\right)  =B.
\]
We can apply this to $A=T_{1}\cap T_{2}$ and $B=T_{1}$ (since $T_{1}\cap
T_{2}$ is a subset of $T_{1}$). We thus obtain
\[
\left(  T_{1}\cap T_{2}\right)  \cap\left(  T_{1}\setminus\left(  T_{1}\cap
T_{2}\right)  \right)  =\varnothing\ \ \ \ \ \ \ \ \ \ \text{and}%
\ \ \ \ \ \ \ \ \ \ \left(  T_{1}\cap T_{2}\right)  \cup\left(  T_{1}%
\setminus\left(  T_{1}\cap T_{2}\right)  \right)  =T_{1}.
\]
In view of $T_{1}\setminus\left(  T_{1}\cap T_{2}\right)  =T_{1}\setminus
T_{2}$, this rewrites as follows:%
\[
\left(  T_{1}\cap T_{2}\right)  \cap\left(  T_{1}\setminus T_{2}\right)
=\varnothing\ \ \ \ \ \ \ \ \ \ \text{and}\ \ \ \ \ \ \ \ \ \ \left(
T_{1}\cap T_{2}\right)  \cup\left(  T_{1}\setminus T_{2}\right)  =T_{1}.
\]

\end{verlong}

\noindent Hence, Theorem \ref{thm.ind.gen-com.split2} (applied to $T_{1}$,
$T_{1}\cap T_{2}$ and $T_{1}\setminus T_{2}$ instead of $S$, $X$ and $Y$)
yields%
\begin{equation}
\sum_{s\in T_{1}}a_{s}=\sum_{s\in T_{1}\cap T_{2}}a_{s}+\underbrace{\sum_{s\in
T_{1}\setminus T_{2}}a_{s}}_{=0}=\sum_{s\in T_{1}\cap T_{2}}a_{s}.
\label{pf.lem.sol.ind.gen-com.fin-sup.proofs.1.6}%
\end{equation}
The same argument (with the roles of $T_{1}$ and $T_{2}$ interchanged) yields%
\[
\sum_{s\in T_{2}}a_{s}=\sum_{s\in T_{2}\cap T_{1}}a_{s}=\sum_{s\in T_{1}\cap
T_{2}}a_{s}%
\]
(since $T_{2}\cap T_{1}=T_{1}\cap T_{2}$). Comparing this equality with
(\ref{pf.lem.sol.ind.gen-com.fin-sup.proofs.1.6}), we obtain $\sum_{s\in
T_{1}}a_{s}=\sum_{s\in T_{2}}a_{s}$. This proves Lemma
\ref{lem.sol.ind.gen-com.fin-sup.proofs.1}.
\end{proof}

\begin{proof}
[Proof of Proposition \ref{prop.ind.gen-com.fin-sup.leg}.]We shall
\textbf{not} use the notation introduced in Definition
\ref{def.ind.gen-com.fin-sup.sum} in this proof, because we have not yet
convinced ourselves that this notation is well-defined.

\textbf{(a)} Lemma \ref{lem.sol.ind.gen-com.fin-sup.proofs.1} shows that if
$T_{1}$ and $T_{2}$ are two finite subsets $T$ of $S$ such that
(\ref{eq.def.ind.gen-com.fin-sup.sum.0}) holds, then
\[
\sum_{s\in T_{1}}a_{s}=\sum_{s\in T_{2}}a_{s}.
\]
In other words, if $T$ is a finite subset of $S$ such that
(\ref{eq.def.ind.gen-com.fin-sup.sum.0}) holds, then the sum $\sum_{s\in
T}a_{s}$ does not depend on the choice of $T$. This proves Proposition
\ref{prop.ind.gen-com.fin-sup.leg} \textbf{(a)}.

\textbf{(b)} Assume that the set $S$ is finite. Note that every $s\in
S\setminus S$ satisfies $a_{s}=0$. (Indeed, this is vacuously true, because
there exists no $s\in S\setminus S$ (since $S\setminus S=\varnothing$).)

Let $T$ be a finite subset of $S$ such that%
\begin{equation}
\text{every }s\in S\setminus T\text{ satisfies }a_{s}=0.
\label{pf.prop.ind.gen-com.fin-sup.leg.b.1}%
\end{equation}
(It is easy to see that such a subset $T$ exists\footnote{\textit{Proof.}
Indeed, $S$ itself is such a subset $T$ (because $S$ is a finite subset of $S$
and has the property that every $s\in S\setminus S$ satisfies $a_{s}=0$).}.)

We shall prove the equality $\sum_{s\in S}a_{s}=\sum_{s\in T}a_{s}$. (Note
that both sums in this equality are defined according to Definition
\ref{def.ind.gen-com.defsum1}, not according to Definition
\ref{def.ind.gen-com.fin-sup.sum}.)

Indeed, we have%
\[
\sum_{s\in S\setminus T}\underbrace{a_{s}}_{\substack{=0\\\text{(by
(\ref{pf.prop.ind.gen-com.fin-sup.leg.b.1}))}}}=\sum_{s\in S\setminus T}0=0
\]
(by Theorem \ref{thm.ind.gen-com.sum(0)} (applied to $S\setminus T$ instead of
$S$)).

Now, $T$ and $S\setminus T$ are two subsets of $S$. A well-known fact from set
theory says that if $A$ is a subset of a set $B$, then
\[
A\cap\left(  B\setminus A\right)  =\varnothing\ \ \ \ \ \ \ \ \ \ \text{and}%
\ \ \ \ \ \ \ \ \ \ A\cup\left(  B\setminus A\right)  =B.
\]
We can apply this to $A=T$ and $B=S$. We thus obtain%
\[
T\cap\left(  S\setminus T\right)  =\varnothing\ \ \ \ \ \ \ \ \ \ \text{and}%
\ \ \ \ \ \ \ \ \ \ T\cup\left(  S\setminus T\right)  =S.
\]
Thus, Theorem \ref{thm.ind.gen-com.split2} (applied to $T$ and $S\setminus T$
instead of $X$ and $Y$) yields%
\begin{equation}
\sum_{s\in S}a_{s}=\sum_{s\in T}a_{s}+\underbrace{\sum_{s\in S\setminus
T}a_{s}}_{=0}=\sum_{s\in T}a_{s}. \label{pf.prop.ind.gen-com.fin-sup.leg.b.4}%
\end{equation}

Now, recall that $T$ is a finite subset of $S$ such that every $s\in
S\setminus T$ satisfies $a_{s}=0$. Hence, Definition
\ref{def.ind.gen-com.fin-sup.sum} defines the sum $\sum_{s\in S}a_{s}$ to be
$\sum_{s\in T}a_{s}$. Thus,%
\begin{align*}
&  \left(  \text{the sum }\sum_{s\in S}a_{s}\text{ defined in Definition
\ref{def.ind.gen-com.fin-sup.sum}}\right) \\
&  =\sum_{s\in T}a_{s}=\sum_{s\in S}a_{s}\ \ \ \ \ \ \ \ \ \ \left(  \text{by
(\ref{pf.prop.ind.gen-com.fin-sup.leg.b.4})}\right) \\
&  =\left(  \text{the sum }\sum_{s\in S}a_{s}\text{ defined in Definition
\ref{def.ind.gen-com.defsum1}}\right)  .
\end{align*}
In other words, the sum $\sum_{s\in S}a_{s}$ defined in Definition
\ref{def.ind.gen-com.fin-sup.sum} is identical with the sum $\sum_{s\in
S}a_{s}$ defined in Definition \ref{def.ind.gen-com.defsum1}. This proves
Proposition \ref{prop.ind.gen-com.fin-sup.leg} \textbf{(b)}.
\end{proof}

Now that Proposition \ref{prop.ind.gen-com.fin-sup.leg} has been proven, we
can use Definition \ref{def.ind.gen-com.fin-sup.sum}, and we shall do so in
the remainder of this section.

Let us state a takeaway from Definition \ref{def.ind.gen-com.fin-sup.sum}:

\begin{corollary}
\label{cor.ind.gen-com.fin-sup.def}Let $S$ be any set. Let $\left(
a_{s}\right)  _{s\in S}$ be a finitely supported $\mathbb{A}$-valued
$S$-family. Let $T$ be a finite subset of $S$ such that%
\[
\text{every }s\in S\setminus T\text{ satisfies }a_{s}=0.
\]
Then,%
\[
\sum_{s\in S}a_{s}=\sum_{s\in T}a_{s}.
\]

\end{corollary}

\begin{proof}
[Proof of Corollary \ref{cor.ind.gen-com.fin-sup.def}.]Recall that $T$ is a
finite subset of $S$ such that every $s\in S\setminus T$ satisfies $a_{s}=0$.
Hence, Definition \ref{def.ind.gen-com.fin-sup.sum} defines the sum
$\sum_{s\in S}a_{s}$ to be $\sum_{s\in T}a_{s}$. Thus, $\sum_{s\in S}%
a_{s}=\sum_{s\in T}a_{s}$. This proves Corollary
\ref{cor.ind.gen-com.fin-sup.def}.
\end{proof}

Next, let us generalize Corollary \ref{cor.ind.gen-com.fin-sup.def}, no longer
requiring $T$ to be finite:

\begin{corollary}
\label{cor.ind.gen-com.fin-sup.remove0}Let $S$ be any set. Let $\left(
a_{s}\right)  _{s\in S}$ be a finitely supported $\mathbb{A}$-valued
$S$-family. Let $T$ be a subset of $S$ such that%
\begin{equation}
\text{every }s\in S\setminus T\text{ satisfies }a_{s}=0.
\label{eq.cor.ind.gen-com.fin-sup.remove0.ass}%
\end{equation}
Then, the $\mathbb{A}$-valued $T$-family $\left(  a_{s}\right)  _{s\in T}$ is
finitely supported as well, and satisfies%
\begin{equation}
\sum_{s\in S}a_{s}=\sum_{s\in T}a_{s}.
\label{eq.cor.ind.gen-com.fin-sup.remove0.claim}%
\end{equation}

\end{corollary}

Note that both sums $\sum_{s\in S}a_{s}$ and $\sum_{s\in T}a_{s}$ in Corollary
\ref{cor.ind.gen-com.fin-sup.remove0} are defined according to Definition
\ref{def.ind.gen-com.fin-sup.sum}.

\begin{proof}
[Proof of Corollary \ref{cor.ind.gen-com.fin-sup.remove0}.]We have assumed
that the $S$-family $\left(  a_{s}\right)  _{s\in S}$ is finitely supported.
In other words, only finitely many $s\in S$ satisfy $a_{s}\neq0$ (by the
definition of \textquotedblleft finitely supported\textquotedblright). In
other words, there exists a finite subset $Q$ of $S$ such that%
\begin{equation}
\text{every }s\in S\setminus Q\text{ satisfies }a_{s}=0.
\label{pf.cor.ind.gen-com.fin-sup.remove0.1}%
\end{equation}
Consider this $Q$.

The set $T\cap Q$ is a subset of $Q$, and thus is finite (since $Q$ is
finite). Also, $T\cap Q$ is a subset of $T$. Moreover, every $s\in
T\setminus\left(  T\cap Q\right)  $ satisfies $a_{s}=0$%
\ \ \ \ \footnote{\textit{Proof.} Let $s\in T\setminus\left(  T\cap Q\right)
$. We must prove that $a_{s}=0$.
\par
We have $s\in T\setminus\left(  T\cap Q\right)  =\underbrace{T}_{\subseteq
S}\setminus Q\subseteq S\setminus Q$. Hence,
(\ref{pf.cor.ind.gen-com.fin-sup.remove0.1}) shows that $a_{s}=0$. Qed.}.
Hence, there exists a finite subset $F$ of $T$ such that every $s\in
T\setminus F$ satisfies $a_{s}=0$ (namely, $F=T\cap Q$). In other words, only
finitely many $s\in T$ satisfy $a_{s}\neq0$. In other words, the $\mathbb{A}%
$-valued $T$-family $\left(  a_{s}\right)  _{s\in T}$ is finitely supported
(by the definition of \textquotedblleft finitely supported\textquotedblright).
It now remains to prove (\ref{eq.cor.ind.gen-com.fin-sup.remove0.claim}).

The set $T\cap Q$ is finite and is a subset of $S$ (since $T\cap Q\subseteq
Q\subseteq S$), and has the property that%
\begin{equation}
\text{every }s\in S\setminus\left(  T\cap Q\right)  \text{ satisfies }a_{s}=0.
\label{pf.cor.ind.gen-com.fin-sup.remove0.2}%
\end{equation}

[\textit{Proof of (\ref{pf.cor.ind.gen-com.fin-sup.remove0.2}):} Let $s\in
S\setminus\left(  T\cap Q\right)  $. We must prove that $a_{s}=0$.

We have $s\in S\setminus\left(  T\cap Q\right)  $. In other words, $s\in S$
and $s\notin T\cap Q$.

We are in one of the following two cases:

\textit{Case 1:} We have $s\in Q$.

\textit{Case 2:} We don't have $s\in Q$.

Let us first consider Case 1. In this case, we have $s\in Q$. Combining $s\in
Q$ with $s\notin T\cap Q$, we obtain $s\in Q\setminus\left(  T\cap Q\right)
=\underbrace{Q}_{\subseteq S}\setminus T\subseteq S\setminus T$. Thus,
$a_{s}=0$ (by (\ref{eq.cor.ind.gen-com.fin-sup.remove0.ass})). Hence,
$a_{s}=0$ is proven in Case 1.

Let us now consider Case 2. In this case, we don't have $s\in Q$. Hence,
$s\notin Q$. Combining $s\in S$ with $s\notin Q$, we obtain $s\in S\setminus
Q$. Hence, (\ref{pf.cor.ind.gen-com.fin-sup.remove0.1}) yields $a_{s}=0$.
Thus, $a_{s}=0$ is proven in Case 2.

We now have shown that $a_{s}=0$ in each of the two Cases 1 and 2. Since these
two Cases cover all possibilities, we thus conclude that $a_{s}=0$ always
holds. This proves (\ref{pf.cor.ind.gen-com.fin-sup.remove0.2}).]

Hence, Corollary \ref{cor.ind.gen-com.fin-sup.def} (applied to $T\cap Q$
instead of $T$) yields%
\begin{equation}
\sum_{s\in S}a_{s}=\sum_{s\in T\cap Q}a_{s}.
\label{pf.cor.ind.gen-com.fin-sup.remove0.L}%
\end{equation}

On the other hand, $T\cap Q$ is a finite subset of $T$, and has the property
that%
\begin{equation}
\text{every }s\in T\setminus\left(  T\cap Q\right)  \text{ satisfies }a_{s}=0.
\label{pf.cor.ind.gen-com.fin-sup.remove0.3}%
\end{equation}

\begin{vershort}
\noindent(Indeed, this follows from
(\ref{pf.cor.ind.gen-com.fin-sup.remove0.1}), because every $s\in
T\setminus\left(  T\cap Q\right)  $ satisfies $s\in T\setminus\left(  T\cap
Q\right)  =\underbrace{T}_{\subseteq S}\setminus Q\subseteq S\setminus Q$.)
\end{vershort}

\begin{verlong}
[\textit{Proof of (\ref{pf.cor.ind.gen-com.fin-sup.remove0.3}):} Let $s\in
T\setminus\left(  T\cap Q\right)  $. We must prove that $a_{s}=0$.

We have $s\in T\setminus\left(  T\cap Q\right)  =\underbrace{T}_{\subseteq
S}\setminus Q\subseteq S\setminus Q$. Hence,
(\ref{pf.cor.ind.gen-com.fin-sup.remove0.1}) yields $a_{s}=0$. Thus,
(\ref{pf.cor.ind.gen-com.fin-sup.remove0.3}) is proven.]
\end{verlong}

Hence, Corollary \ref{cor.ind.gen-com.fin-sup.def} (applied to $T\cap Q$ and
$T$ instead of $T$ and $S$) yields%
\[
\sum_{s\in T}a_{s}=\sum_{s\in T\cap Q}a_{s}%
\]
(since the $\mathbb{A}$-valued $T$-family $\left(  a_{s}\right)  _{s\in T}$ is
finitely supported). Comparing this with
(\ref{pf.cor.ind.gen-com.fin-sup.remove0.L}), we obtain $\sum_{s\in S}%
a_{s}=\sum_{s\in T}a_{s}$. This proves
(\ref{eq.cor.ind.gen-com.fin-sup.remove0.claim}). Hence, Corollary
\ref{cor.ind.gen-com.fin-sup.remove0} is proven.
\end{proof}

\begin{proof}
[Proof of Theorem \ref{thm.ind.gen-com.sum(a+b).gen}.]The $S$-family $\left(
a_{s}\right)  _{s\in S}$ is finitely supported. In other words, only finitely
many $s\in S$ satisfy $a_{s}\neq0$ (by the definition of \textquotedblleft
finitely supported\textquotedblright). In other words, there exists a finite
subset $A$ of $S$ such that%
\begin{equation}
\text{every }s\in S\setminus A\text{ satisfies }a_{s}=0.
\label{pf.thm.ind.gen-com.sum(a+b).gen.A}%
\end{equation}
Consider this $A$.

The $S$-family $\left(  b_{s}\right)  _{s\in S}$ is finitely supported. In
other words, only finitely many $s\in S$ satisfy $b_{s}\neq0$ (by the
definition of \textquotedblleft finitely supported\textquotedblright). In
other words, there exists a finite subset $B$ of $S$ such that%
\begin{equation}
\text{every }s\in S\setminus B\text{ satisfies }b_{s}=0.
\label{pf.thm.ind.gen-com.sum(a+b).gen.B}%
\end{equation}
Consider this $B$.

The set $A\cup B$ is a subset of $S$ (since $A$ and $B$ are subsets of $S$).
Moreover, this set $A\cup B$ is finite (since both sets $A$ and $B$ are
finite). Furthermore,
\begin{equation}
\text{every }s\in S\setminus\left(  A\cup B\right)  \text{ satisfies }a_{s}=0
\label{pf.thm.ind.gen-com.sum(a+b).gen.AuB1}%
\end{equation}
\footnote{\textit{Proof.} Let $s\in S\setminus\left(  A\cup B\right)  $. We
must prove that $a_{s}=0$.
\par
We have $s\in S\setminus\left(  A\cup B\right)  =\left(  S\setminus A\right)
\setminus B\subseteq S\setminus A$. Thus,
(\ref{pf.thm.ind.gen-com.sum(a+b).gen.A}) yields $a_{s}=0$. Qed.}. Also,%
\begin{equation}
\text{every }s\in S\setminus\left(  A\cup B\right)  \text{ satisfies }b_{s}=0
\label{pf.thm.ind.gen-com.sum(a+b).gen.AuB2}%
\end{equation}
\footnote{\textit{Proof.} Let $s\in S\setminus\left(  A\cup B\right)  $. We
must prove that $b_{s}=0$.
\par
We have $s\in S\setminus\left(  A\cup B\right)  =\left(  S\setminus B\right)
\setminus A\subseteq S\setminus B$. Thus,
(\ref{pf.thm.ind.gen-com.sum(a+b).gen.B}) yields $b_{s}=0$. Qed.}. Thus,%
\begin{equation}
\text{every }s\in S\setminus\left(  A\cup B\right)  \text{ satisfies }%
a_{s}+b_{s}=0 \label{pf.thm.ind.gen-com.sum(a+b).gen.AuB3}%
\end{equation}
\footnote{\textit{Proof.} Every $s\in S\setminus\left(  A\cup B\right)  $
satisfies $\underbrace{a_{s}}_{\substack{=0\\\text{(by
(\ref{pf.thm.ind.gen-com.sum(a+b).gen.AuB1}))}}}+\underbrace{b_{s}%
}_{\substack{=0\\\text{(by (\ref{pf.thm.ind.gen-com.sum(a+b).gen.AuB2}))}%
}}=0+0=0$, qed.}. Thus, there exists a finite subset $T$ of $S$ such that
every $s\in S\setminus T$ satisfies $a_{s}+b_{s}=0$ (namely, $T=A\cup B$). In
other words, only finitely many $s\in S$ satisfy $a_{s}+b_{s}\neq0$. In other
words, the $S$-family $\left(  a_{s}+b_{s}\right)  _{s\in S}$ is finitely supported.

The set $A\cup B$ is finite. Hence, Theorem \ref{thm.ind.gen-com.sum(a+b)}
(applied to $A\cup B$ instead of $S$) yields%
\begin{equation}
\sum_{s\in A\cup B}\left(  a_{s}+b_{s}\right)  =\sum_{s\in A\cup B}a_{s}%
+\sum_{s\in A\cup B}b_{s}. \label{pf.thm.ind.gen-com.sum(a+b).gen.finsums}%
\end{equation}

We know that $A\cup B$ is a finite subset of $S$ such that every $s\in
S\setminus\left(  A\cup B\right)  $ satisfies $a_{s}=0$ (by
(\ref{pf.thm.ind.gen-com.sum(a+b).gen.AuB1})). Thus, Corollary
\ref{cor.ind.gen-com.fin-sup.def} (applied to $A\cup B$ instead of $T$) yields%
\begin{equation}
\sum_{s\in S}a_{s}=\sum_{s\in A\cup B}a_{s}.
\label{pf.thm.ind.gen-com.sum(a+b).gen.suma}%
\end{equation}

We know that $A\cup B$ is a finite subset of $S$ such that every $s\in
S\setminus\left(  A\cup B\right)  $ satisfies $b_{s}=0$ (by
(\ref{pf.thm.ind.gen-com.sum(a+b).gen.AuB2})). Thus, Corollary
\ref{cor.ind.gen-com.fin-sup.def} (applied to $\left(  b_{s}\right)  _{s\in
S}$ and $A\cup B$ instead of $\left(  a_{s}\right)  _{s\in S}$ and $T$) yields%
\begin{equation}
\sum_{s\in S}b_{s}=\sum_{s\in A\cup B}b_{s}.
\label{pf.thm.ind.gen-com.sum(a+b).gen.sumb}%
\end{equation}

We know that $A\cup B$ is a finite subset of $S$ such that every $s\in
S\setminus\left(  A\cup B\right)  $ satisfies $a_{s}+b_{s}=0$ (by
(\ref{pf.thm.ind.gen-com.sum(a+b).gen.AuB3})). Thus, Corollary
\ref{cor.ind.gen-com.fin-sup.def} (applied to $\left(  a_{s}+b_{s}\right)
_{s\in S}$ and $A\cup B$ instead of $\left(  a_{s}\right)  _{s\in S}$ and $T$)
yields%
\begin{align*}
\sum_{s\in S}\left(  a_{s}+b_{s}\right)   &  =\sum_{s\in A\cup B}\left(
a_{s}+b_{s}\right)  =\underbrace{\sum_{s\in A\cup B}a_{s}}_{\substack{=\sum
_{s\in S}a_{s}\\\text{(by (\ref{pf.thm.ind.gen-com.sum(a+b).gen.suma}))}%
}}+\underbrace{\sum_{s\in A\cup B}b_{s}}_{\substack{=\sum_{s\in S}%
b_{s}\\\text{(by (\ref{pf.thm.ind.gen-com.sum(a+b).gen.sumb}))}}%
}\ \ \ \ \ \ \ \ \ \ \left(  \text{by
(\ref{pf.thm.ind.gen-com.sum(a+b).gen.finsums})}\right) \\
&  =\sum_{s\in S}a_{s}+\sum_{s\in S}b_{s}.
\end{align*}
This completes the proof of Theorem \ref{thm.ind.gen-com.sum(a+b).gen}.
\end{proof}

Before we prove Theorem \ref{thm.ind.gen-com.sheph}, let us restate Corollary
\ref{cor.ind.gen-com.fin-sup.def} with different names:

\begin{corollary}
\label{cor.ind.gen-com.fin-sup.defW}Let $W$ be any set. Let $\left(
a_{w}\right)  _{w\in W}$ be a finitely supported $\mathbb{A}$-valued
$W$-family. Let $V$ be a finite subset of $W$ such that%
\[
\text{every }w\in W\setminus V\text{ satisfies }a_{w}=0.
\]
Then,%
\[
\sum_{w\in W}a_{w}=\sum_{w\in V}a_{w}.
\]

\end{corollary}

\begin{proof}
[Proof of Corollary \ref{cor.ind.gen-com.fin-sup.defW}.]Corollary
\ref{cor.ind.gen-com.fin-sup.defW} is just Corollary
\ref{cor.ind.gen-com.fin-sup.def}, with all appearances of the letters
\textquotedblleft$S$\textquotedblright, \textquotedblleft$s$\textquotedblright%
\ and \textquotedblleft$T$\textquotedblright\ replaced by \textquotedblleft%
$W$\textquotedblright, \textquotedblleft$w$\textquotedblright\ and
\textquotedblleft$V$\textquotedblright. Thus, Corollary
\ref{cor.ind.gen-com.fin-sup.defW} follows from the latter corollary.
\end{proof}

\begin{proof}
[Proof of Theorem \ref{thm.ind.gen-com.sheph}.]Let $V$ be the subset $f\left(
S\right)  $ of $W$. (Thus, \newline$V=f\left(  S\right)  =\left\{  f\left(
s\right)  \ \mid\ s\in S\right\}  $.) The set $f\left(  S\right)  $ is finite
(since the set $S$ is finite). In other words, the set $V$ is finite (since
$V=f\left(  S\right)  $).

Define a map $g:S\rightarrow V$ by%
\[
\left(  g\left(  s\right)  =f\left(  s\right)  \text{ for each }s\in S\right)
.
\]
(This map $g$ is well-defined, since each $s\in S$ satisfies $f\left(
s\right)  \in f\left(  S\right)  =V$.)

For each $w\in W$, we define an element $b_{w}\in\mathbb{A}$ by%
\begin{equation}
b_{w}=\sum_{\substack{s\in S;\\f\left(  s\right)  =w}}a_{s}.
\label{pf.thm.ind.gen-com.sheph.bf=}%
\end{equation}

Now, each $w\in V$ satisfies%
\begin{equation}
\sum_{\substack{s\in S;\\g\left(  s\right)  =w}}a_{s}=b_{w}.
\label{pf.thm.ind.gen-com.sheph.1}%
\end{equation}

[\textit{Proof of (\ref{pf.thm.ind.gen-com.sheph.1}):} Let $w\in V$. Then,
$g\left(  s\right)  =f\left(  s\right)  $ for each $s\in S$ (by the definition
of $g$). Hence, the summation sign \textquotedblleft$\sum_{\substack{s\in
S;\\g\left(  s\right)  =w}}$\textquotedblright\ can be rewritten as
\textquotedblleft$\sum_{\substack{s\in S;\\f\left(  s\right)  =w}%
}$\textquotedblright. Thus, $\sum_{\substack{s\in S;\\g\left(  s\right)
=w}}a_{s}=\sum_{\substack{s\in S;\\f\left(  s\right)  =w}}a_{s}=b_{w}$ (by
(\ref{pf.thm.ind.gen-com.sheph.bf=})). This proves
(\ref{pf.thm.ind.gen-com.sheph.1}).]

Theorem \ref{thm.ind.gen-com.shephf} (applied to $V$ and $g$ instead of $W$
and $f$) yields%
\begin{equation}
\sum_{s\in S}a_{s}=\sum_{w\in V}\underbrace{\sum_{\substack{s\in S;\\g\left(
s\right)  =w}}a_{s}}_{\substack{=b_{w}\\\text{(by
(\ref{pf.thm.ind.gen-com.sheph.1}))}}}=\sum_{w\in V}b_{w}.
\label{pf.thm.ind.gen-com.sheph.fin}%
\end{equation}

Moreover,%
\begin{equation}
\text{every }w\in W\setminus V\text{ satisfies }b_{w}=0.
\label{pf.thm.ind.gen-com.sheph.2}%
\end{equation}

[\textit{Proof of (\ref{pf.thm.ind.gen-com.sheph.2}):} Let $w\in W\setminus
V$. Thus, $w\in W$ and $w\notin V$. Hence, every $s\in S$ satisfies $f\left(
s\right)  \neq w$\ \ \ \ \footnote{\textit{Proof.} Let $s\in S$. If we had
$f\left(  s\right)  =w$, then we would have $w=f\left(  \underbrace{s}_{\in
S}\right)  \in f\left(  S\right)  =V$, which would contradict $w\notin V$.
Hence, we cannot have $f\left(  s\right)  =w$. In other words, we have
$f\left(  s\right)  \neq w$, qed.}. In other words, there exists no $s\in S$
satisfying $f\left(  s\right)  =w$. Hence, the sum $\sum_{\substack{s\in
S;\\f\left(  s\right)  =w}}a_{s}$ is an empty sum. Thus, $\sum_{\substack{s\in
S;\\f\left(  s\right)  =w}}a_{s}=\left(  \text{empty sum}\right)  =0$. In view
of (\ref{pf.thm.ind.gen-com.sheph.bf=}), this rewrites as $b_{w}=0$. This
proves (\ref{pf.thm.ind.gen-com.sheph.2}).]

Thus, there exists a finite subset $T$ of $W$ such that every $w\in W\setminus
T$ satisfies $b_{w}=0$ (namely, $T=V$). In other words, only finitely many
$w\in W$ satisfy $b_{w}\neq0$. In other words, the $\mathbb{A}$-valued
$W$-family $\left(  b_{w}\right)  _{w\in W}$ is finitely supported. In view of
(\ref{pf.thm.ind.gen-com.sheph.bf=}), this rewrites as follows: The
$\mathbb{A}$-valued $W$-family $\left(  \sum_{\substack{s\in S;\\f\left(
s\right)  =w}}a_{s}\right)  _{w\in W}$ is finitely supported. Hence, the sum
$\sum_{w\in W}\sum_{\substack{s\in S;\\f\left(  s\right)  =w}}a_{s}$ is well-defined.

Moreover, Corollary \ref{cor.ind.gen-com.fin-sup.defW} (applied to $\left(
b_{w}\right)  _{w\in W}$ instead of $\left(  a_{w}\right)  _{w\in W}$) yields
that $\sum_{w\in W}b_{w}=\sum_{w\in V}b_{w}$. Comparing this with
(\ref{pf.thm.ind.gen-com.sheph.fin}), we obtain%
\[
\sum_{s\in S}a_{s}=\sum_{w\in W}\underbrace{b_{w}}_{\substack{=\sum
_{\substack{s\in S;\\f\left(  s\right)  =w}}a_{s}\\\text{(by
(\ref{pf.thm.ind.gen-com.sheph.bf=}))}}}=\sum_{w\in W}\sum_{\substack{s\in
S;\\f\left(  s\right)  =w}}a_{s}.
\]
This completes the proof of Theorem \ref{thm.ind.gen-com.sheph}.
\end{proof}

\begin{proof}
[Proof of Theorem \ref{thm.ind.gen-com.sheph-gen}.]The $S$-family $\left(
a_{s}\right)  _{s\in S}$ is finitely supported. In other words, only finitely
many $s\in S$ satisfy $a_{s}\neq0$. In other words, there exists a finite
subset $T$ of $S$ such that every $s\in S\setminus T$ satisfies $a_{s}=0$.
Consider this $T$, and denote it by $P$. Thus, $P$ is a finite subset of $S$
and has the property that%
\begin{equation}
\text{every }s\in S\setminus P\text{ satisfies }a_{s}=0.
\label{pf.thm.ind.gen-com.sheph-gen.P}%
\end{equation}
Thus, Corollary \ref{cor.ind.gen-com.fin-sup.def} (applied to $T=P$) yields%
\begin{equation}
\sum_{s\in S}a_{s}=\sum_{s\in P}a_{s}.
\label{pf.thm.ind.gen-com.sheph-gen.sum-P}%
\end{equation}

Let $g:P\rightarrow W$ be the restriction $f\mid_{P}$ of the map $f$ to the
subset $P$ of $S$. Hence,
\begin{equation}
g\left(  t\right)  =f\left(  t\right)  \ \ \ \ \ \ \ \ \ \ \text{for each
}t\in P \label{pf.thm.ind.gen-com.sheph-gen.gt=ft}%
\end{equation}
(by the definition of a restriction of a map).

Theorem \ref{thm.ind.gen-com.sheph} (applied to $P$ and $g$ instead of $S$ and
$f$) yields that the $\mathbb{A}$-valued $W$-family $\left(  \sum
_{\substack{s\in P;\\g\left(  s\right)  =w}}a_{s}\right)  _{w\in W}$ is
finitely supported and satisfies%
\begin{equation}
\sum_{s\in P}a_{s}=\sum_{w\in W}\sum_{\substack{s\in P;\\g\left(  s\right)
=w}}a_{s}. \label{pf.thm.ind.gen-com.sheph-gen.fin}%
\end{equation}
Hence, (\ref{pf.thm.ind.gen-com.sheph-gen.sum-P}) becomes%
\begin{equation}
\sum_{s\in S}a_{s}=\sum_{s\in P}a_{s}=\sum_{w\in W}\sum_{\substack{s\in
P;\\g\left(  s\right)  =w}}a_{s}. \label{pf.thm.ind.gen-com.sheph-gen.finR}%
\end{equation}

Now, we claim the following:

\begin{statement}
\textit{Claim 1:} Let $w\in W$.

\textbf{(a)} The $\mathbb{A}$-valued $\left\{  t\in S\ \mid\ f\left(
t\right)  =w\right\}  $-family $\left(  a_{s}\right)  _{s\in\left\{  t\in
S\ \mid\ f\left(  t\right)  =w\right\}  }$ is finitely supported (so that the
sum $\sum_{\substack{s\in S;\\f\left(  s\right)  =w}}a_{s}$ is well-defined).

\textbf{(b)} We have $\sum_{\substack{s\in S;\\f\left(  s\right)  =w}%
}a_{s}=\sum_{\substack{s\in P;\\g\left(  s\right)  =w}}a_{s}$.
\end{statement}

[\textit{Proof of Claim 1:} Define a subset $Q$ of $S$ by $Q=\left\{  t\in
S\ \mid\ f\left(  t\right)  =w\right\}  $.

The set $P\cap Q$ is a subset of $P$, and thus is finite (since $P$ is
finite). Moreover, this set $P\cap Q$ is a subset of $Q$. Finally,
\[
\text{every }s\in Q\setminus\left(  P\cap Q\right)  \text{ satisfies }a_{s}=0
\]
\footnote{\textit{Proof:} Let $s\in Q\setminus\left(  P\cap Q\right)  $. We
must prove that $a_{s}=0$.
\par
We have $s\in Q\setminus\left(  P\cap Q\right)  =\underbrace{Q}_{\subseteq
S}\setminus P\subseteq S\setminus P$. Hence,
(\ref{pf.thm.ind.gen-com.sheph-gen.P}) shows that $a_{s}=0$. Qed.}. Hence,
there exists a finite subset $T$ of $Q$ such that every $s\in Q\setminus T$
satisfies $a_{s}=0$ (namely, $T=P\cap Q$). In other words, only finitely many
$s\in Q$ satisfy $a_{s}\neq0$. In other words, the $\mathbb{A}$-valued
$Q$-family $\left(  a_{s}\right)  _{s\in Q}$ is finitely supported. In view of
$Q=\left\{  t\in S\ \mid\ f\left(  t\right)  =w\right\}  $, this rewrites as
follows: The $\mathbb{A}$-valued $\left\{  t\in S\ \mid\ f\left(  t\right)
=w\right\}  $-family $\left(  a_{s}\right)  _{s\in\left\{  t\in S\ \mid
\ f\left(  t\right)  =w\right\}  }$ is finitely supported. This proves Claim 1
\textbf{(a)}.

\textbf{(b)} We know that $P\cap Q$ is a subset of $Q$ such that every $s\in
Q\setminus\left(  P\cap Q\right)  $ satisfies $a_{s}=0$. Hence, Corollary
\ref{cor.ind.gen-com.fin-sup.remove0} (applied to $Q$ and $P\cap Q$ instead of
$S$ and $T$) shows that the $\mathbb{A}$-valued $P\cap Q$-family $\left(
a_{s}\right)  _{s\in P\cap Q}$ is finitely supported as well, and satisfies%
\begin{equation}
\sum_{s\in Q}a_{s}=\sum_{s\in P\cap Q}a_{s}.
\label{pf.thm.ind.gen-com.sheph-gen.c1.pf.1}%
\end{equation}

But $Q=\left\{  t\in S\ \mid\ f\left(  t\right)  =w\right\}  $. Hence, an
element $t$ of $S$ belongs to $Q$ if and only if it satisfies $f\left(
t\right)  =w$. In other words, for any element $t\in S$, we have the logical
equivalence%
\begin{equation}
\left(  t\in Q\right)  \ \Longleftrightarrow\ \left(  f\left(  t\right)
=w\right)  . \label{pf.thm.ind.gen-com.sheph-gen.c1.pf.2}%
\end{equation}
Now, the definition of the intersection of two sets shows that
\begin{align*}
P\cap Q  &  =\left\{  t\in P\ \mid\ \underbrace{t\in Q}%
_{\substack{\Longleftrightarrow\ \left(  f\left(  t\right)  =w\right)
\\\text{(by (\ref{pf.thm.ind.gen-com.sheph-gen.c1.pf.2})}\\\text{(since }t\in
P\subseteq S\text{))}}}\right\}  =\left\{  t\in P\ \mid\ \underbrace{f\left(
t\right)  }_{\substack{=g\left(  t\right)  \\\text{(by
(\ref{pf.thm.ind.gen-com.sheph-gen.gt=ft}))}}}=w\right\} \\
&  =\left\{  t\in P\ \mid\ g\left(  t\right)  =w\right\}  .
\end{align*}
Now,%
\begin{align*}
\underbrace{\sum_{\substack{s\in S;\\f\left(  s\right)  =w}}}_{\substack{=\sum
_{s\in\left\{  t\in S\ \mid\ f\left(  t\right)  =w\right\}  }=\sum_{s\in
Q}\\\text{(since }\left\{  t\in S\ \mid\ f\left(  t\right)  =w\right\}
=Q\text{)}}}a_{s}  &  =\sum_{s\in Q}a_{s}=\underbrace{\sum_{s\in P\cap Q}%
}_{=\sum_{\substack{s\in\left\{  t\in P\ \mid\ g\left(  t\right)  =w\right\}
\\\text{(since }P\cap Q=\left\{  t\in P\ \mid\ g\left(  t\right)  =w\right\}
\text{)}}}}a_{s}\ \ \ \ \ \ \ \ \ \ \left(  \text{by
(\ref{pf.thm.ind.gen-com.sheph-gen.c1.pf.1})}\right) \\
&  =\underbrace{\sum_{s\in\left\{  t\in P\ \mid\ g\left(  t\right)
=w\right\}  }}_{=\sum_{\substack{s\in P;\\g\left(  s\right)  =w}}}a_{s}%
=\sum_{\substack{s\in P;\\g\left(  s\right)  =w}}a_{s}.
\end{align*}
This proves Claim 1 \textbf{(b)}.]

Now, Claim 1 \textbf{(a)} shows that for each $w\in W$, the $\mathbb{A}%
$-valued $\left\{  t\in S\ \mid\ f\left(  t\right)  =w\right\}  $-family
$\left(  a_{s}\right)  _{s\in\left\{  t\in S\ \mid\ f\left(  t\right)
=w\right\}  }$ is finitely supported (so that the sum $\sum_{\substack{s\in
S;\\f\left(  s\right)  =w}}a_{s}$ is well-defined). Also, Claim 1 \textbf{(b)}
shows that%
\[
\sum_{\substack{s\in S;\\f\left(  s\right)  =w}}a_{s}=\sum_{\substack{s\in
P;\\g\left(  s\right)  =w}}a_{s}\ \ \ \ \ \ \ \ \ \ \text{for each }w\in W.
\]
In other words,%
\[
\left(  \sum_{\substack{s\in S;\\f\left(  s\right)  =w}}a_{s}\right)  _{w\in
W}=\left(  \sum_{\substack{s\in P;\\g\left(  s\right)  =w}}a_{s}\right)
_{w\in W}.
\]
Hence, the $\mathbb{A}$-valued $W$-family $\left(  \sum_{\substack{s\in
S;\\f\left(  s\right)  =w}}a_{s}\right)  _{w\in W}$ is finitely supported
(since we know that the $\mathbb{A}$-valued $W$-family $\left(  \sum
_{\substack{s\in P;\\g\left(  s\right)  =w}}a_{s}\right)  _{w\in W}$ is
finitely supported). Finally, (\ref{pf.thm.ind.gen-com.sheph-gen.finR})
becomes%
\[
\sum_{s\in S}a_{s}=\sum_{w\in W}\underbrace{\sum_{\substack{s\in P;\\g\left(
s\right)  =w}}a_{s}}_{\substack{=\sum_{\substack{s\in S;\\f\left(  s\right)
=w}}a_{s}\\\text{(by Claim 1 \textbf{(b)})}}}=\sum_{w\in W}\sum
_{\substack{s\in S;\\f\left(  s\right)  =w}}a_{s}.
\]
This completes the proof of Theorem \ref{thm.ind.gen-com.sheph-gen}.
\end{proof}

\begin{proof}
[Solution to Exercise \ref{exe.ind.gen-com.fin-sup.proofs}.]We have proven
Proposition \ref{prop.ind.gen-com.fin-sup.leg}, Theorem
\ref{thm.ind.gen-com.sum(a+b).gen}, Theorem \ref{thm.ind.gen-com.sheph} and
Theorem \ref{thm.ind.gen-com.sheph-gen}. Thus, Exercise
\ref{exe.ind.gen-com.fin-sup.proofs} is solved.
\end{proof}

\subsection{Solution to Exercise \ref{exe.mod.parity}}

\begin{proof}
[Proof of Proposition \ref{prop.mod.even-odd-mod2}.]\textbf{(a)} We have the
following chain of logical equivalences:%
\begin{align*}
&  \left(  \text{the integer }n\text{ is even}\right) \\
\Longleftrightarrow\  &  \left(  n\text{ is divisible by }2\right)
\ \ \ \ \ \ \ \ \ \ \left(  \text{by the definition of \textquotedblleft
even\textquotedblright}\right) \\
\Longleftrightarrow\  &  \left(  2\mid n\right)  \ \Longleftrightarrow
\ \left(  2\mid n-0\right)  \ \ \ \ \ \ \ \ \ \ \left(  \text{since
}n=n-0\right) \\
\Longleftrightarrow\  &  \left(  n\equiv0\operatorname{mod}2\right)
\ \ \ \ \ \ \ \ \ \ \left(  \text{by the definition of \textquotedblleft%
}n\equiv0\operatorname{mod}2\text{\textquotedblright}\right)  .
\end{align*}
In other words, the integer $n$ is even if and only if $n\equiv
0\operatorname{mod}2$. This proves Proposition \ref{prop.mod.even-odd-mod2}
\textbf{(a)}.

\textbf{(b)} Proposition \ref{prop.ind.quo-rem.odd} \textbf{(a)} shows that
the integer $n$ is odd if and only if $n$ can be written in the form $n=2m+1$
for some $m\in\mathbb{Z}$.

Now, we claim the following logical implication:%
\begin{equation}
\left(  \text{the integer }n\text{ is odd}\right)  \ \Longrightarrow\ \left(
n\equiv1\operatorname{mod}2\right)  . \label{pf.prop.mod.even-odd-mod2.b.1}%
\end{equation}

[\textit{Proof of (\ref{pf.prop.mod.even-odd-mod2.b.1}):} Assume that the
integer $n$ is odd. We shall show that $n\equiv1\operatorname{mod}2$.

Recall that the integer $n$ is odd if and only if $n$ can be written in the
form $n=2m+1$ for some $m\in\mathbb{Z}$. Hence, $n$ can be written in the form
$n=2m+1$ for some $m\in\mathbb{Z}$ (since the integer $n$ is odd). Consider
this $m$. Thus, $n=2m+1$, so that $n-1=2m$. Thus, $n-1$ is divisible by $2$.
In other words, $n\equiv1\operatorname{mod}2$ (by the definition of
\textquotedblleft$n\equiv1\operatorname{mod}2$\textquotedblright). This proves
the implication (\ref{pf.prop.mod.even-odd-mod2.b.1}).]

Next, we claim the following logical implication:%
\begin{equation}
\left(  n\equiv1\operatorname{mod}2\right)  \ \Longrightarrow\ \left(
\text{the integer }n\text{ is odd}\right)  .
\label{pf.prop.mod.even-odd-mod2.b.2}%
\end{equation}

[\textit{Proof of (\ref{pf.prop.mod.even-odd-mod2.b.2}):} Assume that
$n\equiv1\operatorname{mod}2$. We shall show that the integer $n$ is odd.

We have $n\equiv1\operatorname{mod}2$. In other words, $n-1$ is divisible by
$2$ (by the definition of \textquotedblleft$n\equiv1\operatorname{mod}%
2$\textquotedblright). In other words, there exists a $z\in\mathbb{Z}$ such
that $n-1=2z$. Consider this $z$. Now, $n-1=2z$, so that $n=2z+1$. Hence, $n$
can be written in the form $n=2m+1$ for some $m\in\mathbb{Z}$ (namely, for
$m=z$). Hence, the integer $n$ is odd (since the integer $n$ is odd if and
only if $n$ can be written in the form $n=2m+1$ for some $m\in\mathbb{Z}$).
This proves the implication (\ref{pf.prop.mod.even-odd-mod2.b.2}).]

Combining the two implications (\ref{pf.prop.mod.even-odd-mod2.b.1}) and
(\ref{pf.prop.mod.even-odd-mod2.b.2}), we obtain the equivalence
\[
\left(  \text{the integer }n\text{ is odd}\right)  \ \Longleftrightarrow
\ \left(  n\equiv1\operatorname{mod}2\right)  .
\]
In other words, the integer $n$ is odd if and only if $n\equiv
1\operatorname{mod}2$. This proves Proposition \ref{prop.mod.even-odd-mod2}
\textbf{(b)}.
\end{proof}

\begin{proof}
[Proof of Proposition \ref{prop.mod.parity}.]Let $u\%2$ denote the remainder
of the division of $u$ by $2$. Then, Corollary \ref{cor.ind.quo-rem.remmod}
\textbf{(a)} (applied to $N=2$ and $n=u$) yields that $u\%2\in\left\{
0,1,\ldots,2-1\right\}  $ and $u\%2\equiv u\operatorname{mod}2$. Thus,
$u\%2\in\left\{  0,1,\ldots,2-1\right\}  =\left\{  0,1\right\}  $.

Let $v\%2$ denote the remainder of the division of $v$ by $2$. Then, Corollary
\ref{cor.ind.quo-rem.remmod} \textbf{(a)} (applied to $N=2$ and $n=v$) yields
that $v\%2\in\left\{  0,1,\ldots,2-1\right\}  $ and $v\%2\equiv
v\operatorname{mod}2$. Thus, $v\%2\in\left\{  0,1,\ldots,2-1\right\}
=\left\{  0,1\right\}  $.

We have either $u\%2=0$ or $u\%2=1$ (since $u\%2\in\left\{  0,1\right\}  $).
Thus, we are in one of the following two cases:

\textit{Case 1:} We have $u\%2=0$.

\textit{Case 2:} We have $u\%2=1$.

Let us first consider Case 1. In this case, we have $u\%2=0$.

From $u\%2\equiv u\operatorname{mod}2$, we obtain $u\equiv
u\%2=0\operatorname{mod}2$. But Proposition \ref{prop.mod.even-odd-mod2}
\textbf{(a)} (applied to $n=u$) shows that the integer $u$ is even if and only
if $u\equiv0\operatorname{mod}2$. Thus, $u$ is even (since $u\equiv
0\operatorname{mod}2$). Hence, Corollary \ref{cor.mod.-1powers} \textbf{(a)}
(applied to $n=u$) yields $\left(  -1\right)  ^{u}=1$.

We have either $v\%2=0$ or $v\%2=1$ (since $v\%2\in\left\{  0,1\right\}  $).
Thus, we are in one of the following two subcases:

\textit{Subcase 1.1:} We have $v\%2=0$.

\textit{Subcase 1.2:} We have $v\%2=1$.

Let us first consider Subcase 1.1. In this case, we have $v\%2=0$. But
$u\equiv0=v\%2\equiv v\operatorname{mod}2$. Thus, the statement $\left(
u\equiv v\operatorname{mod}2\right)  $ is true. Also, $v\equiv
v\%2=0\operatorname{mod}2$.

Proposition \ref{prop.mod.even-odd-mod2} \textbf{(a)} (applied to $n=v$) shows
that the integer $v$ is even if and only if $v\equiv0\operatorname{mod}2$.
Thus, $v$ is even (since $v\equiv0\operatorname{mod}2$). Hence, Corollary
\ref{cor.mod.-1powers} \textbf{(a)} (applied to $n=v$) yields $\left(
-1\right)  ^{v}=1$. Thus, $\left(  -1\right)  ^{u}=1=\left(  -1\right)  ^{v}$.
Hence, the statement $\left(  \left(  -1\right)  ^{u}=\left(  -1\right)
^{v}\right)  $ is true.

Finally, the numbers $u$ and $v$ are either both even or both odd (since $u$
and $v$ are both even). Thus, the statement%
\[
\left(  u\text{ and }v\text{ are either both even or both odd}\right)
\]
is true.

We have now shown that the three statements
\[
\left(  u\equiv v\operatorname{mod}2\right)  ,\ \left(  u\text{ and }v\text{
are either both even or both odd}\right)  \text{ and }\left(  \left(
-1\right)  ^{u}=\left(  -1\right)  ^{v}\right)
\]
are all true. Thus, these three statements are equivalent. In other words, we
have
\begin{align*}
\left(  u\equiv v\operatorname{mod}2\right)  \  &  \Longleftrightarrow
\ \left(  u\text{ and }v\text{ are either both even or both odd}\right) \\
&  \Longleftrightarrow\ \left(  \left(  -1\right)  ^{u}=\left(  -1\right)
^{v}\right)  .
\end{align*}
Hence, Proposition \ref{prop.mod.parity} is proven in Subcase 1.1.

Let us next consider Subcase 1.2. In this case, we have $v\%2=1$. From
$v\%2\equiv v\operatorname{mod}2$, we obtain $v\equiv v\%2=1\operatorname{mod}%
2$. From $u\equiv0\operatorname{mod}2$, we obtain $0\equiv u\operatorname{mod}%
2$. Hence, if we had $u\equiv v\operatorname{mod}2$, then we would have
$0\equiv u\equiv v\equiv1\operatorname{mod}2$, which would contradict the fact
that $0\not \equiv 1\operatorname{mod}2$ (since $2\nmid0-1$). Thus, we cannot
have $u\equiv v\operatorname{mod}2$. In other words, the statement $\left(
u\equiv v\operatorname{mod}2\right)  $ is false.

Proposition \ref{prop.mod.even-odd-mod2} \textbf{(b)} (applied to $n=v$) shows
that the integer $v$ is odd if and only if $v\equiv1\operatorname{mod}2$.
Thus, $v$ is odd (since $v\equiv1\operatorname{mod}2$). Hence, Corollary
\ref{cor.mod.-1powers} \textbf{(b)} (applied to $n=v$) yields $\left(
-1\right)  ^{v}=-1$. Thus, $\left(  -1\right)  ^{u}=1\neq-1=\left(  -1\right)
^{v}$. Hence, the statement $\left(  \left(  -1\right)  ^{u}=\left(
-1\right)  ^{v}\right)  $ is false.

Finally, the number $v$ is not even (since $v$ is odd), while the number $u$
is not odd (since $u$ is even). Hence, the numbers $u$ and $v$ are neither
both even nor both odd. Thus, the statement%
\[
\left(  u\text{ and }v\text{ are either both even or both odd}\right)
\]
is false.

We have now shown that the three statements
\[
\left(  u\equiv v\operatorname{mod}2\right)  ,\ \left(  u\text{ and }v\text{
are either both even or both odd}\right)  \text{ and }\left(  \left(
-1\right)  ^{u}=\left(  -1\right)  ^{v}\right)
\]
are all false. Thus, these three statements are equivalent. In other words, we
have
\begin{align*}
\left(  u\equiv v\operatorname{mod}2\right)  \  &  \Longleftrightarrow
\ \left(  u\text{ and }v\text{ are either both even or both odd}\right) \\
&  \Longleftrightarrow\ \left(  \left(  -1\right)  ^{u}=\left(  -1\right)
^{v}\right)  .
\end{align*}
Hence, Proposition \ref{prop.mod.parity} is proven in Subcase 1.2.

We have now proven Proposition \ref{prop.mod.parity} in each of the two
Subcases 1.1 and 1.2. Since these two Subcases cover the whole Case 1, we thus
conclude that Proposition \ref{prop.mod.parity} holds in Case 1.

\begin{vershort}
The arguments needed to deal with Case 2 are very similar: We again must split
the case into two subcases, the first being $v\%2=1$ and the second being
$v\%2=0$; the reasoning in each subcase is completely analogous to the
reasoning in the corresponding subcase of Case 1 (but with the roles of
\textquotedblleft even\textquotedblright\ and \textquotedblleft
odd\textquotedblright\ interchanged). We leave the details to the reader.
\end{vershort}

\begin{verlong}
Let us next consider Case 2. In this case, we have $u\%2=1$.

From $u\%2\equiv u\operatorname{mod}2$, we obtain $u\equiv
u\%2=1\operatorname{mod}2$. But Proposition \ref{prop.mod.even-odd-mod2}
\textbf{(b)} (applied to $n=u$) shows that the integer $u$ is odd if and only
if $u\equiv1\operatorname{mod}2$. Thus, $u$ is odd (since $u\equiv
1\operatorname{mod}2$). Hence, Corollary \ref{cor.mod.-1powers} \textbf{(b)}
(applied to $n=u$) yields $\left(  -1\right)  ^{u}=-1$.

We have either $v\%2=1$ or $v\%2=0$ (since $v\%2\in\left\{  0,1\right\}
=\left\{  1,0\right\}  $). Thus, we are in one of the following two subcases:

\textit{Subcase 2.1:} We have $v\%2=1$.

\textit{Subcase 2.2:} We have $v\%2=0$.

Let us first consider Subcase 2.1. In this case, we have $v\%2=1$. From
$u\equiv1\operatorname{mod}2$, we obtain $u\equiv1=v\%2\equiv
v\operatorname{mod}2$. Thus, the statement $\left(  u\equiv
v\operatorname{mod}2\right)  $ is true. Also, $v\equiv
v\%2=1\operatorname{mod}2$.

Proposition \ref{prop.mod.even-odd-mod2} \textbf{(b)} (applied to $n=v$) shows
that the integer $v$ is odd if and only if $v\equiv1\operatorname{mod}2$.
Thus, $v$ is odd (since $v\equiv1\operatorname{mod}2$). Hence, Corollary
\ref{cor.mod.-1powers} \textbf{(b)} (applied to $n=v$) yields $\left(
-1\right)  ^{v}=-1$. Thus, $\left(  -1\right)  ^{u}=-1=\left(  -1\right)
^{v}$. Hence, the statement $\left(  \left(  -1\right)  ^{u}=\left(
-1\right)  ^{v}\right)  $ is true.

Finally, the numbers $u$ and $v$ are either both even or both odd (since $u$
and $v$ are both odd). Thus, the statement%
\[
\left(  u\text{ and }v\text{ are either both even or both odd}\right)
\]
is true.

We have now shown that the three statements
\[
\left(  u\equiv v\operatorname{mod}2\right)  ,\ \left(  u\text{ and }v\text{
are either both even or both odd}\right)  \text{ and }\left(  \left(
-1\right)  ^{u}=\left(  -1\right)  ^{v}\right)
\]
are all true. Thus, these three statements are equivalent. In other words, we
have
\begin{align*}
\left(  u\equiv v\operatorname{mod}2\right)  \  &  \Longleftrightarrow
\ \left(  u\text{ and }v\text{ are either both even or both odd}\right) \\
&  \Longleftrightarrow\ \left(  \left(  -1\right)  ^{u}=\left(  -1\right)
^{v}\right)  .
\end{align*}
Hence, Proposition \ref{prop.mod.parity} is proven in Subcase 2.1.

Let us next consider Subcase 2.2. In this case, we have $v\%2=0$. From
$v\%2\equiv v\operatorname{mod}2$, we obtain $v\equiv v\%2=0\operatorname{mod}%
2$. From $u\equiv1\operatorname{mod}2$, we obtain $1\equiv u\operatorname{mod}%
2$. Hence, if we had $u\equiv v\operatorname{mod}2$, then we would have
$1\equiv u\equiv v\equiv0\operatorname{mod}2$, which would contradict the fact
that $1\not \equiv 0\operatorname{mod}2$ (since $2\nmid1-0$). Thus, we cannot
have $u\equiv v\operatorname{mod}2$. In other words, the statement $\left(
u\equiv v\operatorname{mod}2\right)  $ is false.

Proposition \ref{prop.mod.even-odd-mod2} \textbf{(a)} (applied to $n=v$) shows
that the integer $v$ is even if and only if $v\equiv0\operatorname{mod}2$.
Thus, $v$ is even (since $v\equiv0\operatorname{mod}2$). Hence, Corollary
\ref{cor.mod.-1powers} \textbf{(a)} (applied to $n=v$) yields $\left(
-1\right)  ^{v}=1$. Thus, $\left(  -1\right)  ^{u}=-1\neq1=\left(  -1\right)
^{v}$. Hence, the statement $\left(  \left(  -1\right)  ^{u}=\left(
-1\right)  ^{v}\right)  $ is false.

Finally, the number $u$ is not even (since $u$ is odd), while the number $v$
is not odd (since $v$ is even). Hence, the numbers $u$ and $v$ are neither
both even nor both odd. Thus, the statement%
\[
\left(  u\text{ and }v\text{ are either both even or both odd}\right)
\]
is false.

We have now shown that the three statements
\[
\left(  u\equiv v\operatorname{mod}2\right)  ,\ \left(  u\text{ and }v\text{
are either both even or both odd}\right)  \text{ and }\left(  \left(
-1\right)  ^{u}=\left(  -1\right)  ^{v}\right)
\]
are all false. Thus, these three statements are equivalent. In other words, we
have
\begin{align*}
\left(  u\equiv v\operatorname{mod}2\right)  \  &  \Longleftrightarrow
\ \left(  u\text{ and }v\text{ are either both even or both odd}\right) \\
&  \Longleftrightarrow\ \left(  \left(  -1\right)  ^{u}=\left(  -1\right)
^{v}\right)  .
\end{align*}
Hence, Proposition \ref{prop.mod.parity} is proven in Subcase 2.2.

We have now proven Proposition \ref{prop.mod.parity} in each of the two
Subcases 2.1 and 2.2. Since these two Subcases cover the whole Case 2, we thus
conclude that Proposition \ref{prop.mod.parity} holds in Case 2.
\end{verlong}

We have now proven Proposition \ref{prop.mod.parity} in each of the two Cases
1 and 2. Since these two Cases cover all possibilities, we thus conclude that
Proposition \ref{prop.mod.parity} always holds.
\end{proof}

\begin{proof}
[Solution to Exercise \ref{exe.mod.parity}.]We have proven Proposition
\ref{prop.mod.parity} and Proposition \ref{prop.mod.even-odd-mod2}; thus,
Exercise \ref{exe.mod.parity} is solved.
\end{proof}

\subsection{Solution to Exercise \ref{exe.mod.unique-cong}}

\begin{proof}
[Solution to Exercise \ref{exe.mod.unique-cong}.]Let $n=h-p-1$. Thus,
$n\in\mathbb{Z}$. Let $n\%N$ denote the remainder of the division of $n$ by
$N$. Thus, Corollary \ref{cor.ind.quo-rem.remmod} \textbf{(a)} yields that
$n\%N\in\left\{  0,1,\ldots,N-1\right\}  $ and $n\%N\equiv n\operatorname{mod}%
N$.

Let $x=\left(  p+1\right)  +\left(  n\%N\right)  $. Then, $x\in\mathbb{Z}$
(since both $p+1$ and $n\%N$ belong to $\mathbb{Z}$) and
\begin{align*}
x  &  =\left(  p+1\right)  +\left(  n\%N\right) \\
&  \in\left\{  \left(  p+1\right)  +0,\left(  p+1\right)  +1,\ldots,\left(
p+1\right)  +\left(  N-1\right)  \right\} \\
&  \ \ \ \ \ \ \ \ \ \ \left(  \text{since }n\%N\in\left\{  0,1,\ldots
,N-1\right\}  \right) \\
&  =\left\{  p+1,p+2,\ldots,p+N\right\}  .
\end{align*}
Moreover, $x=\left(  p+1\right)  +\underbrace{\left(  n\%N\right)  }_{\equiv
n=h-p-1\operatorname{mod}N}\equiv\left(  p+1\right)  +\left(  h-p-1\right)
=h\operatorname{mod}N$. Hence, $x$ is an element $g\in\left\{  p+1,p+2,\ldots
,p+N\right\}  $ satisfying $g\equiv h\operatorname{mod}N$ (since $x\in\left\{
p+1,p+2,\ldots,p+N\right\}  $ and $x\equiv h\operatorname{mod}N$). Thus, there
exists \textbf{at least one} element $g\in\left\{  p+1,p+2,\ldots,p+N\right\}
$ satisfying $g\equiv h\operatorname{mod}N$ (namely, $g=x$).

Now, let $g\in\left\{  p+1,p+2,\ldots,p+N\right\}  $ be such that $g\equiv
h\operatorname{mod}N$. We shall prove that $g=x$.

Indeed, $\underbrace{g}_{\equiv h\operatorname{mod}N}-\left(  p+1\right)
\equiv h-\left(  p+1\right)  =h-p-1=n\operatorname{mod}N$. Also, from
$g\in\left\{  p+1,p+2,\ldots,p+N\right\}  $, we obtain $g-\left(  p+1\right)
\in\left\{  0,1,\ldots,N-1\right\}  $. Thus, Corollary
\ref{cor.ind.quo-rem.remmod} \textbf{(c)} (applied to $c=g-\left(  p+1\right)
$) yields $g-\left(  p+1\right)  =n\%N$. Hence, $g=\left(  p+1\right)
+\left(  n\%N\right)  =x$.

Now, forget that we fixed $g$. We thus have shown that if $g\in\left\{
p+1,p+2,\ldots,p+N\right\}  $ satisfies $g\equiv h\operatorname{mod}N$, then
$g=x$. In other words, every element \newline$g\in\left\{  p+1,p+2,\ldots
,p+N\right\}  $ satisfying $g\equiv h\operatorname{mod}N$ must be equal to
$x$. Thus, there exists \textbf{at most one} element $g\in\left\{
p+1,p+2,\ldots,p+N\right\}  $ satisfying $g\equiv h\operatorname{mod}N$.

We have now shown the following two facts:

\begin{itemize}
\item There exists \textbf{at least one} element $g\in\left\{  p+1,p+2,\ldots
,p+N\right\}  $ satisfying $g\equiv h\operatorname{mod}N$.

\item There exists \textbf{at most one} element $g\in\left\{  p+1,p+2,\ldots
,p+N\right\}  $ satisfying $g\equiv h\operatorname{mod}N$.
\end{itemize}

Combining these two facts, we conclude that there exists a \textbf{unique}
element $g\in\left\{  p+1,p+2,\ldots,p+N\right\}  $ satisfying $g\equiv
h\operatorname{mod}N$. This solves Exercise \ref{exe.mod.unique-cong}.
\end{proof}

\subsection{Solution to Exercise \ref{exe.mod.ak-bk}}

\begin{proof}
[Solution to Exercise \ref{exe.mod.ak-bk}.]\textbf{(a)} We have $a-b=1\left(
a-b\right)  $, thus $a-b\mid a-b$. In other words, $a\equiv
b\operatorname{mod}a-b$ (by the definition of \textquotedblleft
congruent\textquotedblright). Hence, Proposition \ref{prop.mod.pow} (applied
to $n=a-b$) yields $a^{k}\equiv b^{k}\operatorname{mod}a-b$. In other words,
$a-b\mid a^{k}-b^{k}$ (by the definition of \textquotedblleft
congruent\textquotedblright). This solves Exercise \ref{exe.mod.ak-bk}
\textbf{(a)}.

\textbf{(b)} We know that $k$ is odd. Hence, Corollary \ref{cor.mod.-1powers}
\textbf{(b)} (applied to $n=k$) yields $\left(  -1\right)  ^{k}=-1$.

Exercise \ref{exe.mod.ak-bk} \textbf{(a)} (applied to $-b$ instead of $b$)
yields $a-\left(  -b\right)  \mid a^{k}-\left(  -b\right)  ^{k}$. In view of
$a-\left(  -b\right)  =a+b$ and%
\[
a^{k}-\underbrace{\left(  -b\right)  ^{k}}_{=\left(  -1\right)  ^{k}b^{k}%
}=a^{k}-\underbrace{\left(  -1\right)  ^{k}}_{=-1}b^{k}=a^{k}-\left(
-1\right)  b^{k}=a^{k}+b^{k},
\]
this rewrites as $a+b\mid a^{k}+b^{k}$. This solves Exercise
\ref{exe.mod.ak-bk} \textbf{(b)}.
\end{proof}

\subsection{Solution to Exercise \ref{exe.ind.LP2.another-div}}

\begin{proof}
[Solution to Exercise \ref{exe.ind.LP2.another-div}.]Proposition
\ref{prop.ind.LP2} \textbf{(a)} shows that we have%
\begin{equation}
b_{n}\in\mathbb{N}\text{ for each }n\in\mathbb{N}\text{.}
\label{sol.ind.LP2.another-div.inN}%
\end{equation}

For every integer $m\geq1$, we have%
\begin{equation}
b_{m}^{r}+1=b_{m+1}b_{m-1}. \label{sol.ind.LP2.another-div.o0}%
\end{equation}
(Indeed, this is precisely the identity (\ref{pf.prop.ind.LP2.o0}), which was
already proven during our proof of Proposition \ref{prop.ind.LP2}.)

\begin{vershort}
Now, let $n$ be a positive integer. Thus, the four integers $n-1,n,n+1,n+2$
all belong to $\mathbb{N}$. Hence, the numbers $b_{n-1},b_{n},b_{n+1},b_{n+2}$
all belong to $\mathbb{N}$ (by (\ref{sol.ind.LP2.another-div.inN})).
\end{vershort}

\begin{verlong}
Now, let $n$ be a positive integer. Thus, $n-1\in\mathbb{N}$. Hence,
(\ref{sol.ind.LP2.another-div.inN}) (applied to $n-1$ instead of $n$) yields
$b_{n-1}\in\mathbb{N}$. Also, (\ref{sol.ind.LP2.another-div.inN}) (applied to
$n+2$ instead of $\mathbb{N}$) yields $b_{n+2}\in\mathbb{N}$. Also,
(\ref{sol.ind.LP2.another-div.inN}) (applied to $n+1$ instead of $n$) yields
$b_{n+1}\in\mathbb{N}$. Finally, (\ref{sol.ind.LP2.another-div.inN}) yields
$b_{n}\in\mathbb{N}$. Hence, $b_{n}\in\mathbb{N}\subseteq\mathbb{Z}$.

Thus, all four numbers $b_{n-1},b_{n},b_{n+1},b_{n+2}$ belong to $\mathbb{N}$,
and therefore are integers. This allows us to state congruences that involve
these four numbers.
\end{verlong}

We have $n\geq1$ (since $n$ is a positive integer). Hence,
(\ref{sol.ind.LP2.another-div.o0}) (applied to $m=n$) yields $b_{n}%
^{r}+1=b_{n+1}b_{n-1}$.

Define an integer $x$ by $x=b_{n}+1$. (This $x$ is indeed an integer, because
$b_{n}\in\mathbb{N}\subseteq\mathbb{Z}$.) We have\footnote{We recall that the
four numbers $b_{n-1},b_{n},b_{n+1},b_{n+2}$ are integers; this justifies the
use of these four numbers in congruences modulo $x$.} $b_{n}\equiv
-1\operatorname{mod}x$ (since $b_{n}-\left(  -1\right)  =b_{n}+1=x$ is
divisible by $x$). Also, Exercise \ref{exe.mod.ak-bk} \textbf{(b)} (applied to
$r$, $b_{n}$ and $1$ instead of $k$, $a$ and $b$) shows that $b_{n}+1\mid
b_{n}^{r}+1$. In view of $b_{n}+1=x$ and $b_{n}^{r}+1=b_{n+1}b_{n-1}$, this
rewrites as $x\mid b_{n+1}b_{n-1}$. In other words, $b_{n+1}b_{n-1}%
\equiv0\operatorname{mod}x$. But $r\neq0$ (since $r$ is odd). Combining this
with $r\geq0$, we obtain $r>0$. Hence, $r\geq1$ (since $r$ is an integer).
Now,%
\begin{equation}
b_{n-1}\underbrace{b_{n+1}^{r}}_{\substack{=b_{n+1}b_{n+1}^{r-1}\\\text{(since
}r\geq1\text{)}}}=\underbrace{b_{n-1}b_{n+1}}_{\equiv0\operatorname{mod}%
x}b_{n+1}^{r-1}\equiv0\operatorname{mod}x. \label{sol.ind.LP2.another-div.2}%
\end{equation}

Now, $n+1\geq n\geq1$. Thus, (\ref{sol.ind.LP2.another-div.o0}) (applied to
$m=n+1$) yields
\begin{equation}
b_{n+1}^{r}+1=\underbrace{b_{\left(  n+1\right)  +1}}_{=b_{n+2}}%
\underbrace{b_{\left(  n+1\right)  -1}}_{=b_{n}}=b_{n+2}b_{n}.
\label{sol.ind.LP2.another-div.3}%
\end{equation}

We have%
\begin{align*}
b_{n}b_{n-1}\left(  b_{n+2}+1\right)   &  =b_{n-1}\underbrace{b_{n+2}b_{n}%
}_{\substack{=b_{n+1}^{r}+1\\\text{(by (\ref{sol.ind.LP2.another-div.3}))}%
}}+\underbrace{b_{n}}_{\equiv-1\operatorname{mod}x}b_{n-1}\\
&  \equiv b_{n-1}\left(  b_{n+1}^{r}+1\right)  +\left(  -1\right)
b_{n-1}=b_{n-1}b_{n+1}^{r}\equiv0\operatorname{mod}x
\end{align*}
(by (\ref{sol.ind.LP2.another-div.2})). But%
\begin{align*}
\underbrace{x}_{=b_{n}+1}b_{n-1}\left(  b_{n+2}+1\right)   &  =\left(
b_{n}+1\right)  b_{n-1}\left(  b_{n-2}+1\right) \\
&  =\underbrace{b_{n}b_{n-1}\left(  b_{n+2}+1\right)  }_{\equiv
0\operatorname{mod}x}+b_{n-1}\left(  b_{n-2}+1\right)  \equiv b_{n-1}\left(
b_{n-2}+1\right)  \operatorname{mod}x.
\end{align*}
Hence,%
\[
b_{n-1}\left(  b_{n-2}+1\right)  \equiv\underbrace{x}_{\equiv
0\operatorname{mod}x}b_{n-1}\left(  b_{n+2}+1\right)  \equiv
0\operatorname{mod}x.
\]
In other words, $x\mid b_{n-1}\left(  b_{n-2}+1\right)  $. In view of
$x=b_{n}+1$, this rewrites as $b_{n}+1\mid b_{n-1}\left(  b_{n+2}+1\right)  $.
This solves Exercise \ref{exe.ind.LP2.another-div}.
\end{proof}

\subsection{Solution to Exercise \ref{exe.ind.backw}}

To prove Theorem \ref{thm.ind.IPg-}, we just mildly adapt our proof of Theorem
\ref{thm.ind.IPg}:

\begin{proof}
[Proof of Theorem \ref{thm.ind.IPg-}.]For any $n\in\mathbb{N}$, we have
$g-n\in\mathbb{Z}_{\leq g}$\ \ \ \ \footnote{\textit{Proof.} Let
$n\in\mathbb{N}$. Thus, $n\geq0$, so that $g-\underbrace{n}_{\geq0}\leq
g-0=g$. Hence, $g-n$ is an integer $\leq g$. In other words, $g-n\in
\mathbb{Z}_{\leq g}$ (since $\mathbb{Z}_{\leq g}$ is the set of all integers
that are $\leq g$). Qed.}. Hence, for each $n\in\mathbb{N}$, we can define a
logical statement $\mathcal{B}\left(  n\right)  $ by%
\[
\mathcal{B}\left(  n\right)  =\mathcal{A}\left(  g-n\right)  .
\]
Consider this $\mathcal{B}\left(  n\right)  $.

Now, let us consider the Assumptions A and B from Corollary
\ref{cor.ind.IP0.renamed}. We claim that both of these assumptions are satisfied.

Indeed, the statement $\mathcal{A}\left(  g\right)  $ holds (by Assumption 1).
But the definition of the statement $\mathcal{B}\left(  0\right)  $ shows that
$\mathcal{B}\left(  0\right)  =\mathcal{A}\left(  g-0\right)  =\mathcal{A}%
\left(  g\right)  $. Hence, the statement $\mathcal{B}\left(  0\right)  $
holds (since the statement $\mathcal{A}\left(  g\right)  $ holds). In other
words, Assumption A is satisfied.

Now, we shall show that Assumption B is satisfied. Indeed, let $p\in
\mathbb{N}$ be such that $\mathcal{B}\left(  p\right)  $ holds. The definition
of the statement $\mathcal{B}\left(  p\right)  $ shows that $\mathcal{B}%
\left(  p\right)  =\mathcal{A}\left(  g-p\right)  $. Hence, the statement
$\mathcal{A}\left(  g-p\right)  $ holds (since $\mathcal{B}\left(  p\right)  $ holds).

Also, $p\in\mathbb{N}$, so that $p\geq0$ and thus $g-p\leq g$. In other words,
$g-p\in\mathbb{Z}_{\leq g}$ (since $\mathbb{Z}_{\leq g}$ is the set of all
integers that are $\leq g$).

Recall that Assumption 2 holds. In other words, if $m\in\mathbb{Z}_{\leq g}$
is such that $\mathcal{A}\left(  m\right)  $ holds, then $\mathcal{A}\left(
m-1\right)  $ also holds. Applying this to $m=g-p$, we conclude that
$\mathcal{A}\left(  \left(  g-p\right)  -1\right)  $ holds (since
$\mathcal{A}\left(  g-p\right)  $ holds).

But the definition of $\mathcal{B}\left(  p+1\right)  $ yields $\mathcal{B}%
\left(  p+1\right)  =\mathcal{A}\left(  \underbrace{g-\left(  p+1\right)
}_{=\left(  g-p\right)  -1}\right)  =\mathcal{A}\left(  \left(  g-p\right)
-1\right)  $. Hence, the statement $\mathcal{B}\left(  p+1\right)  $ holds
(since the statement $\mathcal{A}\left(  \left(  g-p\right)  -1\right)  $ holds).

Now, forget that we fixed $p$. We thus have shown that if $p\in\mathbb{N}$ is
such that $\mathcal{B}\left(  p\right)  $ holds, then $\mathcal{B}\left(
p+1\right)  $ also holds. In other words, Assumption B is satisfied.

We now know that both Assumption A and Assumption B are satisfied. Hence,
Corollary \ref{cor.ind.IP0.renamed} shows that
\begin{equation}
\mathcal{B}\left(  n\right)  \text{ holds for each }n\in\mathbb{N}.
\label{pf.thm.ind.IPg-.at}%
\end{equation}

Now, let $n\in\mathbb{Z}_{\leq g}$. Thus, $n$ is an integer such that $n\leq
g$ (by the definition of $\mathbb{Z}_{\leq g}$). Hence, $g-n\geq0$, so that
$g-n\in\mathbb{N}$. Thus, (\ref{pf.thm.ind.IPg-.at}) (applied to $g-n$ instead
of $n$) yields that $\mathcal{B}\left(  g-n\right)  $ holds. But the
definition of $\mathcal{B}\left(  g-n\right)  $ yields $\mathcal{B}\left(
g-n\right)  =\mathcal{A}\left(  \underbrace{g-\left(  g-n\right)  }%
_{=n}\right)  =\mathcal{A}\left(  n\right)  $. Hence, the statement
$\mathcal{A}\left(  n\right)  $ holds (since $\mathcal{B}\left(  g-n\right)  $ holds).

Now, forget that we fixed $n$. We thus have shown that $\mathcal{A}\left(
n\right)  $ holds for each $n\in\mathbb{Z}_{\leq g}$. This proves Theorem
\ref{thm.ind.IPg-}.
\end{proof}

Let us restate Theorem \ref{thm.ind.IPg-} as follows:

\begin{corollary}
\label{cor.ind.IPg-.renamed}Let $h\in\mathbb{Z}$. Let $\mathbb{Z}_{\leq h}$ be
the set $\left\{  h,h-1,h-2,\ldots\right\}  $ (that is, the set of all
integers that are $\leq h$). For each $n\in\mathbb{Z}_{\geq h}$, let
$\mathcal{B}\left(  n\right)  $ be a logical statement.

Assume the following:

\begin{statement}
\textit{Assumption A:} The statement $\mathcal{B}\left(  h\right)  $ holds.
\end{statement}

\begin{statement}
\textit{Assumption B:} If $p\in\mathbb{Z}_{\leq h}$ is such that
$\mathcal{B}\left(  p\right)  $ holds, then $\mathcal{B}\left(  p-1\right)  $
also holds.
\end{statement}

Then, $\mathcal{B}\left(  n\right)  $ holds for each $n\in\mathbb{Z}_{\leq h}$.
\end{corollary}

\begin{proof}
[Proof of Corollary \ref{cor.ind.IPg-.renamed}.]Corollary
\ref{cor.ind.IPg-.renamed} is exactly Theorem \ref{thm.ind.IPg-}, except that
some names have been changed:

\begin{itemize}
\item The variable $g$ has been renamed as $h$.

\item The statements $\mathcal{A}\left(  n\right)  $ have been renamed as
$\mathcal{B}\left(  n\right)  $.

\item Assumption 1 and Assumption 2 have been renamed as Assumption A and
Assumption B.

\item The variable $m$ in Assumption B has been renamed as $p$.
\end{itemize}

Thus, Corollary \ref{cor.ind.IPg-.renamed} holds (since Theorem
\ref{thm.ind.IPg-} holds).
\end{proof}

In order to prove Theorem \ref{thm.ind.IPgh-}, we modify our proof of Theorem
\ref{thm.ind.IPgh} as follows:

\begin{proof}
[Proof of Theorem \ref{thm.ind.IPgh-}.]Let $\mathbb{Z}_{\leq h}$ be the set
$\left\{  h,h-1,h-2,\ldots\right\}  $ (that is, the set of all integers that
are $\leq h$).

For each $n\in\mathbb{Z}_{\leq h}$, we define $\mathcal{B}\left(  n\right)  $
to be the logical statement%
\[
\left(  \text{if }n\in\left\{  g,g+1,\ldots,h\right\}  \text{, then
}\mathcal{A}\left(  n\right)  \text{ holds}\right)  .
\]

Now, let us consider the Assumptions A and B from Corollary
\ref{cor.ind.IPg-.renamed}. We claim that both of these assumptions are satisfied.

Assumption 1 says that if $g\leq h$, then the statement $\mathcal{A}\left(
h\right)  $ holds. Thus, $\mathcal{B}\left(  h\right)  $
holds\footnote{\textit{Proof.} Assume that $h\in\left\{  g,g+1,\ldots
,h\right\}  $. Thus, $g\leq h$. But Assumption 1 says that if $g\leq h$, then
the statement $\mathcal{A}\left(  h\right)  $ holds. Hence, the statement
$\mathcal{A}\left(  h\right)  $ holds (since $g\leq h$).
\par
Now, forget that we assumed that $h\in\left\{  g,g+1,\ldots,h\right\}  $. We
thus have proven that if $h\in\left\{  g,g+1,\ldots,h\right\}  $, then
$\mathcal{A}\left(  h\right)  $ holds. In other words, $\mathcal{B}\left(
h\right)  $ holds (because the statement $\mathcal{B}\left(  h\right)  $ is
defined as $\left(  \text{if }h\in\left\{  g,g+1,\ldots,h\right\}  \text{,
then }\mathcal{A}\left(  h\right)  \text{ holds}\right)  $). Qed.}. In other
words, Assumption A is satisfied.

Next, we shall prove that Assumption B is satisfied. Indeed, let
$p\in\mathbb{Z}_{\leq h}$ be such that $\mathcal{B}\left(  p\right)  $ holds.
We shall now show that $\mathcal{B}\left(  p-1\right)  $ also holds.

Indeed, assume that $p-1\in\left\{  g,g+1,\ldots,h\right\}  $. Thus, $p-1\geq
g$, so that $p\geq p-1\geq g$. Combining this with $p\leq h$ (since
$p\in\mathbb{Z}_{\leq h}$), we conclude that $p\in\left\{  g,g+1,\ldots
,h\right\}  $ (since $p$ is an integer). But we have assumed that
$\mathcal{B}\left(  p\right)  $ holds. In other words,%
\[
\text{if }p\in\left\{  g,g+1,\ldots,h\right\}  \text{, then }\mathcal{A}%
\left(  p\right)  \text{ holds}%
\]
(because the statement $\mathcal{B}\left(  p\right)  $ is defined as $\left(
\text{if }p\in\left\{  g,g+1,\ldots,h\right\}  \text{, then }\mathcal{A}%
\left(  p\right)  \text{ holds}\right)  $). Thus, $\mathcal{A}\left(
p\right)  $ holds (since we have $p\in\left\{  g,g+1,\ldots,h\right\}  $).
Also, from $p-1\geq g$, we obtain $p\geq g+1$. Combining this with $p\leq h$,
we find $p\in\left\{  g+1,g+2,\ldots,h\right\}  $. Thus, we know that
$p\in\left\{  g+1,g+2,\ldots,h\right\}  $ is such that $\mathcal{A}\left(
p\right)  $ holds. Hence, Assumption 2 (applied to $m=p$) shows that
$\mathcal{A}\left(  p-1\right)  $ also holds.

Now, forget that we assumed that $p-1\in\left\{  g,g+1,\ldots,h\right\}  $. We
thus have proven that if $p-1\in\left\{  g,g+1,\ldots,h\right\}  $, then
$\mathcal{A}\left(  p-1\right)  $ holds. In other words, $\mathcal{B}\left(
p-1\right)  $ holds (since the statement $\mathcal{B}\left(  p-1\right)  $ is
defined as \newline$\left(  \text{if }p-1\in\left\{  g,g+1,\ldots,h\right\}
\text{, then }\mathcal{A}\left(  p-1\right)  \text{ holds}\right)  $).

Now, forget that we fixed $p$. We thus have proven that if $p\in
\mathbb{Z}_{\leq h}$ is such that $\mathcal{B}\left(  p\right)  $ holds, then
$\mathcal{B}\left(  p-1\right)  $ also holds. In other words, Assumption B is satisfied.

We now know that both Assumption A and Assumption B are satisfied. Hence,
Corollary \ref{cor.ind.IPg-.renamed} shows that
\begin{equation}
\mathcal{B}\left(  n\right)  \text{ holds for each }n\in\mathbb{Z}_{\leq h}.
\label{pf.thm.ind.IPgh-.at}%
\end{equation}

Now, let $n\in\left\{  g,g+1,\ldots,h\right\}  $. Thus, $n\leq h$, so that
$n\in\mathbb{Z}_{\leq h}$. Hence, (\ref{pf.thm.ind.IPgh-.at}) shows that
$\mathcal{B}\left(  n\right)  $ holds. In other words,
\[
\text{if }n\in\left\{  g,g+1,\ldots,h\right\}  \text{, then }\mathcal{A}%
\left(  n\right)  \text{ holds}%
\]
(since the statement $\mathcal{B}\left(  n\right)  $ was defined as $\left(
\text{if }n\in\left\{  g,g+1,\ldots,h\right\}  \text{, then }\mathcal{A}%
\left(  n\right)  \text{ holds}\right)  $). Thus, $\mathcal{A}\left(
n\right)  $ holds (since we have $n\in\left\{  g,g+1,\ldots,h\right\}  $).

Now, forget that we fixed $n$. We thus have shown that $\mathcal{A}\left(
n\right)  $ holds for each $n\in\left\{  g,g+1,\ldots,h\right\}  $. This
proves Theorem \ref{thm.ind.IPgh-}.
\end{proof}

\begin{proof}
[Solution to Exercise \ref{exe.ind.backw}.]Theorem \ref{thm.ind.IPg-} and
Theorem \ref{thm.ind.IPgh-} have been proven. Thus, Exercise
\ref{exe.ind.backw} is solved.
\end{proof}

\subsection{Solution to Exercise \ref{exe.sum--2choosek-cases}}

\begin{proof}
[First solution to Exercise \ref{exe.sum--2choosek-cases}.]We shall solve
Exercise \ref{exe.sum--2choosek-cases} by induction on $n$:

\textit{Induction base:} Comparing
\begin{align*}%
\begin{cases}
0/2+1, & \text{if }0\text{ is even};\\
-\left(  0+1\right)  /2, & \text{if }0\text{ is odd}%
\end{cases}
&  =0/2+1\ \ \ \ \ \ \ \ \ \ \left(  \text{since }0\text{ is even}\right) \\
&  =1
\end{align*}
with $\sum_{k=0}^{0}\left(  -1\right)  ^{k}\left(  k+1\right)
=\underbrace{\left(  -1\right)  ^{0}}_{=1}\underbrace{\left(  0+1\right)
}_{=1}=1$, we find
\[
\sum_{k=0}^{0}\left(  -1\right)  ^{k}\left(  k+1\right)  =%
\begin{cases}
0/2+1, & \text{if }0\text{ is even};\\
-\left(  0+1\right)  /2, & \text{if }0\text{ is odd}%
\end{cases}
.
\]
In other words, Exercise \ref{exe.sum--2choosek-cases} holds for $n=0$. This
completes the induction base.

\textit{Induction step:} Let $m\in\mathbb{Z}_{\geq1}$. Assume that Exercise
\ref{exe.sum--2choosek-cases} holds for $n=m-1$. We must prove that Exercise
\ref{exe.sum--2choosek-cases} holds for $n=m$.

We have assumed that Exercise \ref{exe.sum--2choosek-cases} holds for $n=m-1$.
In other words,%
\begin{equation}
\sum_{k=0}^{m-1}\left(  -1\right)  ^{k}\left(  k+1\right)  =%
\begin{cases}
\left(  m-1\right)  /2+1, & \text{if }m-1\text{ is even};\\
-\left(  \left(  m-1\right)  +1\right)  /2, & \text{if }m-1\text{ is odd}%
\end{cases}
. \label{sol.exe.sum--2choosek-cases.1st.IH}%
\end{equation}

We must prove that Exercise \ref{exe.sum--2choosek-cases} holds for $n=m$. In
other words, we must prove that%
\begin{equation}
\sum_{k=0}^{m}\left(  -1\right)  ^{k}\left(  k+1\right)  =%
\begin{cases}
m/2+1, & \text{if }m\text{ is even};\\
-\left(  m+1\right)  /2, & \text{if }m\text{ is odd}%
\end{cases}
. \label{sol.exe.sum--2choosek-cases.1st.IG}%
\end{equation}

The integer $m$ is either even or odd. Thus, we are in one of the following
two cases:

\textit{Case 1:} The integer $m$ is even.

\textit{Case 2:} The integer $m$ is odd.

Let us first consider Case 1. In this case, the integer $m$ is even. Thus, the
integer $m-1$ is odd\footnote{\textit{Proof.} Proposition
\ref{prop.mod.even-odd-mod2} \textbf{(a)} (applied to $n=m$) shows that the
integer $m$ is even if and only if $m\equiv0\operatorname{mod}2$. Hence,
$m\equiv0\operatorname{mod}2$ (since the integer $m$ is even). Thus,
$\underbrace{m}_{\equiv0\operatorname{mod}2}-1\equiv0-1=-1\equiv
1\operatorname{mod}2$ (since $2\mid\left(  -1\right)  -1$).
\par
But Proposition \ref{prop.mod.even-odd-mod2} \textbf{(b)} (applied to $n=m-1$)
shows that the integer $m-1$ is odd if and only if $m-1\equiv
1\operatorname{mod}2$. Hence, the integer $m-1$ is odd (since $m-1\equiv
1\operatorname{mod}2$).}. Hence, (\ref{sol.exe.sum--2choosek-cases.1st.IH})
becomes%
\begin{align*}
\sum_{k=0}^{m-1}\left(  -1\right)  ^{k}\left(  k+1\right)   &  =%
\begin{cases}
\left(  m-1\right)  /2+1, & \text{if }m-1\text{ is even};\\
-\left(  \left(  m-1\right)  +1\right)  /2, & \text{if }m-1\text{ is odd}%
\end{cases}
\\
&  =-\left(  \left(  m-1\right)  +1\right)  /2\ \ \ \ \ \ \ \ \ \ \left(
\text{since }m-1\text{ is odd}\right)  .
\end{align*}
Also, Corollary \ref{cor.mod.-1powers} \textbf{(a)} (applied to $n=m$) yields
$\left(  -1\right)  ^{m}=1$ (since $m$ is even).

Now, we can split off the addend for $k=m$ from the sum $\sum_{k=0}^{m}\left(
-1\right)  ^{k}\left(  k+1\right)  $ (since $m\in\left\{  0,1,\ldots
,m\right\}  $). We thus obtain%
\begin{align}
\sum_{k=0}^{m}\left(  -1\right)  ^{k}\left(  k+1\right)   &  =\underbrace{\sum
_{k=0}^{m-1}\left(  -1\right)  ^{k}\left(  k+1\right)  }_{=-\left(  \left(
m-1\right)  +1\right)  /2}+\underbrace{\left(  -1\right)  ^{m}}_{=1}\left(
m+1\right) \nonumber\\
&  =-\left(  \left(  m-1\right)  +1\right)  /2+\left(  m+1\right)
=m/2+1.\nonumber
\end{align}
Comparing this with%
\[%
\begin{cases}
m/2+1, & \text{if }m\text{ is even};\\
-\left(  m+1\right)  /2, & \text{if }m\text{ is odd}%
\end{cases}
=m/2+1\ \ \ \ \ \ \ \ \ \ \left(  \text{since }m\text{ is even}\right)  ,
\]
we obtain%
\[
\sum_{k=0}^{m}\left(  -1\right)  ^{k}\left(  k+1\right)  =%
\begin{cases}
m/2+1, & \text{if }m\text{ is even};\\
-\left(  m+1\right)  /2, & \text{if }m\text{ is odd}%
\end{cases}
.
\]
Thus, (\ref{sol.exe.sum--2choosek-cases.1st.IG}) holds in Case 1.

Let us next consider Case 1. In this case, the integer $m$ is odd. Thus, the
integer $m-1$ is even\footnote{\textit{Proof.} Proposition
\ref{prop.mod.even-odd-mod2} \textbf{(b)} (applied to $n=m$) shows that the
integer $m$ is odd if and only if $m\equiv1\operatorname{mod}2$. Hence,
$m\equiv1\operatorname{mod}2$ (since the integer $m$ is odd). Thus,
$\underbrace{m}_{\equiv1\operatorname{mod}2}-1\equiv1-1=0\operatorname{mod}2$.
\par
But Proposition \ref{prop.mod.even-odd-mod2} \textbf{(a)} (applied to $n=m-1$)
shows that the integer $m-1$ is even if and only if $m-1\equiv
0\operatorname{mod}2$. Hence, the integer $m-1$ is even (since $m-1\equiv
0\operatorname{mod}2$).}. Hence, (\ref{sol.exe.sum--2choosek-cases.1st.IH})
becomes%
\begin{align*}
\sum_{k=0}^{m-1}\left(  -1\right)  ^{k}\left(  k+1\right)   &  =%
\begin{cases}
\left(  m-1\right)  /2+1, & \text{if }m-1\text{ is even};\\
-\left(  \left(  m-1\right)  +1\right)  /2, & \text{if }m-1\text{ is odd}%
\end{cases}
\\
&  =\left(  m-1\right)  /2+1\ \ \ \ \ \ \ \ \ \ \left(  \text{since }m-1\text{
is even}\right)  .
\end{align*}
Also, Corollary \ref{cor.mod.-1powers} \textbf{(b)} (applied to $n=m$) yields
$\left(  -1\right)  ^{m}=-1$ (since $m$ is odd).

Now, we can split off the addend for $k=m$ from the sum $\sum_{k=0}^{m}\left(
-1\right)  ^{k}\left(  k+1\right)  $ (since $m\in\left\{  0,1,\ldots
,m\right\}  $). We thus obtain%
\begin{align}
\sum_{k=0}^{m}\left(  -1\right)  ^{k}\left(  k+1\right)   &  =\underbrace{\sum
_{k=0}^{m-1}\left(  -1\right)  ^{k}\left(  k+1\right)  }_{=\left(  m-1\right)
/2+1}+\underbrace{\left(  -1\right)  ^{m}}_{=-1}\left(  m+1\right) \nonumber\\
&  =\left(  m-1\right)  /2+1+\left(  -1\right)  \left(  m+1\right)  =-\left(
m+1\right)  /2.\nonumber
\end{align}
Comparing this with%
\[%
\begin{cases}
m/2+1, & \text{if }m\text{ is even};\\
-\left(  m+1\right)  /2, & \text{if }m\text{ is odd}%
\end{cases}
=-\left(  m+1\right)  /2\ \ \ \ \ \ \ \ \ \ \left(  \text{since }m\text{ is
odd}\right)  ,
\]
we obtain%
\[
\sum_{k=0}^{m}\left(  -1\right)  ^{k}\left(  k+1\right)  =%
\begin{cases}
m/2+1, & \text{if }m\text{ is even};\\
-\left(  m+1\right)  /2, & \text{if }m\text{ is odd}%
\end{cases}
.
\]
Thus, (\ref{sol.exe.sum--2choosek-cases.1st.IG}) holds in Case 2.

We thus have proven that (\ref{sol.exe.sum--2choosek-cases.1st.IG}) holds in
each of the two Cases 1 and 2. Thus, (\ref{sol.exe.sum--2choosek-cases.1st.IG}%
) always holds (since Cases 1 and 2 cover all possibilities). In other words,
Exercise \ref{exe.sum--2choosek-cases} holds for $n=m$. This completes the
induction step. Thus, Exercise \ref{exe.sum--2choosek-cases} is proven by induction.
\end{proof}

\begin{proof}
[Second solution to Exercise \ref{exe.sum--2choosek-cases}.]We have%
\begin{align}
&  \sum_{k=0}^{n}\left(  -1\right)  ^{k}\left(  k+1\right) \nonumber\\
&  =\left(  -1\right)  ^{0}1+\left(  -1\right)  ^{1}2+\left(  -1\right)
^{2}3+\left(  -1\right)  ^{3}4+\cdots+\left(  -1\right)  ^{n}\left(
n+1\right) \nonumber\\
&  =1-2+3-4\pm\cdots+\left(  -1\right)  ^{n}\left(  n+1\right)  .
\label{sol.sum--2choosek-cases.2nd.1}%
\end{align}
Notice that the sign of the last addend on the right hand side depends on
whether $n$ is even or odd. Thus, we distinguish between the following two cases:

\textit{Case 1:} The integer $n$ is even.

\textit{Case 2:} The integer $n$ is odd.

Let us first consider Case 1. In this case, the integer $n$ is even. Thus,
$\left(  -1\right)  ^{n}=1$ (by Corollary \ref{cor.mod.-1powers}
\textbf{(a)}). Hence, (\ref{sol.sum--2choosek-cases.2nd.1}) becomes%
\begin{align*}
\sum_{k=0}^{n}\left(  -1\right)  ^{k}\left(  k+1\right)   &  =1-2+3-4\pm
\cdots+\underbrace{\left(  -1\right)  ^{n}}_{=1}\left(  n+1\right) \\
&  =1-2+3-4\pm\cdots+\left(  n+1\right) \\
&  =\underbrace{\left(  1-2\right)  }_{=-1}+\underbrace{\left(  3-4\right)
}_{=-1}+\underbrace{\left(  5-6\right)  }_{=-1}+\cdots+\underbrace{\left(
\left(  n-1\right)  -n\right)  }_{=-1}+\left(  n+1\right) \\
&  =\underbrace{\underbrace{\left(  \left(  -1\right)  +\left(  -1\right)
+\left(  -1\right)  +\cdots+\left(  -1\right)  \right)  }_{n/2\text{ addends}%
}}_{=n/2\cdot\left(  -1\right)  =-n/2}+\left(  n+1\right) \\
&  =-n/2+\left(  n+1\right)  =n/2+1.
\end{align*}
Comparing this with%
\[%
\begin{cases}
n/2+1, & \text{if }n\text{ is even};\\
-\left(  n+1\right)  /2, & \text{if }n\text{ is odd}%
\end{cases}
=n/2+1\ \ \ \ \ \ \ \ \ \ \left(  \text{since }n\text{ is even}\right)  ,
\]
we obtain $\sum_{k=0}^{n}\left(  -1\right)  ^{k}\left(  k+1\right)  =%
\begin{cases}
n/2+1, & \text{if }n\text{ is even};\\
-\left(  n+1\right)  /2, & \text{if }n\text{ is odd}%
\end{cases}
$. Hence, Exercise \ref{exe.sum--2choosek-cases} is solved in Case 1.

Let us next consider Case 2. In this case, the integer $n$ is odd. Thus,
$\left(  -1\right)  ^{n}=-1$ (by Corollary \ref{cor.mod.-1powers}
\textbf{(b)}). Hence, (\ref{sol.sum--2choosek-cases.2nd.1}) becomes%
\begin{align*}
\sum_{k=0}^{n}\left(  -1\right)  ^{k}\left(  k+1\right)   &  =1-2+3-4\pm
\cdots+\underbrace{\left(  -1\right)  ^{n}}_{=-1}\left(  n+1\right) \\
&  =1-2+3-4\pm\cdots-\left(  n+1\right) \\
&  =\underbrace{\left(  1-2\right)  }_{=-1}+\underbrace{\left(  3-4\right)
}_{=-1}+\underbrace{\left(  5-6\right)  }_{=-1}+\cdots+\underbrace{\left(
n-\left(  n+1\right)  \right)  }_{=-1}\\
&  =\underbrace{\left(  \left(  -1\right)  +\left(  -1\right)  +\left(
-1\right)  +\cdots+\left(  -1\right)  \right)  }_{\left(  n+1\right)  /2\text{
addends}}\\
&  =\left(  n+1\right)  /2\cdot\left(  -1\right)  =-\left(  n+1\right)  /2.
\end{align*}
Comparing this with%
\[%
\begin{cases}
n/2+1, & \text{if }n\text{ is even};\\
-\left(  n+1\right)  /2, & \text{if }n\text{ is odd}%
\end{cases}
=-\left(  n+1\right)  /2\ \ \ \ \ \ \ \ \ \ \left(  \text{since }n\text{ is
odd}\right)  ,
\]
we obtain $\sum_{k=0}^{n}\left(  -1\right)  ^{k}\left(  k+1\right)  =%
\begin{cases}
n/2+1, & \text{if }n\text{ is even};\\
-\left(  n+1\right)  /2, & \text{if }n\text{ is odd}%
\end{cases}
$. Hence, Exercise \ref{exe.sum--2choosek-cases} is solved in Case 2.

We thus have solved Exercise \ref{exe.sum--2choosek-cases} in both Cases 1 and
2. Hence, Exercise \ref{exe.sum--2choosek-cases} always holds.
\end{proof}

\subsection{Solution to Exercise \ref{exe.multinom1}}

\begin{vershort}
\begin{proof}
[Solution to Exercise \ref{exe.multinom1}.]For each $i\in\left\{
0,1,\ldots,m\right\}  $, define an integer $s_{i}\in\mathbb{N}$ by
$s_{i}=k_{1}+k_{2}+\cdots+k_{i}$. We shall prove that%
\begin{equation}
\dfrac{\left(  k_{1}+k_{2}+\cdots+k_{m}\right)  !}{k_{1}!k_{2}!\cdots k_{m}%
!}=\prod_{i=1}^{m}\dbinom{s_{i}}{k_{i}}. \label{sol.multinom1.short.goal}%
\end{equation}

[\textit{Proof of (\ref{sol.multinom1.short.goal}):} If $m=0$, then
(\ref{sol.multinom1.short.goal}) holds\footnote{\textit{Proof.} This is just
an exercise in understanding empty sums and empty products: Empty sums such as
$k_{1}+k_{2}+\cdots+k_{0}$ are defined to be $0$, and empty products such as
$k_{1}!k_{2}!\cdots k_{0}!$ or $\prod_{i=1}^{0}\dbinom{s_{i}}{k_{i}}$ are
defined to be $1$. With this in mind (and remembering that $0!=1$), it becomes
completely straightforward to verify (\ref{sol.multinom1.short.goal}) when
$m=0$.}. Hence, for the rest of the proof of (\ref{sol.multinom1.short.goal}),
we WLOG assume that we don't have $m=0$.

We have $m\geq1$ (since $m\in\mathbb{N}$ but we don't have $m=0$), so that
$m-1\geq0$.

The definition of $s_{0}$ yields $s_{0}=k_{1}+k_{2}+\cdots+k_{0}=\left(
\text{empty sum}\right)  =0$. Hence, $s_{0}!=0!=1$. Now,%
\begin{align}
\prod_{i=1}^{m}\dfrac{1}{s_{i-1}!}  &  =\prod_{i=0}^{m-1}\dfrac{1}{s_{i}%
!}\ \ \ \ \ \ \ \ \ \ \left(  \text{here, we have substituted }i\text{ for
}i-1\text{ in the product}\right) \nonumber\\
&  =\underbrace{\dfrac{1}{s_{0}!}}_{\substack{=1\\\text{(since }%
s_{0}!=1\text{)}}}\prod_{i=1}^{m-1}\dfrac{1}{s_{i}!}%
\ \ \ \ \ \ \ \ \ \ \left(  \text{since }m-1\geq0\right) \nonumber\\
&  =\prod_{i=1}^{m-1}\dfrac{1}{s_{i}!}=\dfrac{1}{\prod_{i=1}^{m-1}s_{i}!}.
\label{sol.multinom1.short.goal.pf.2}%
\end{align}

The definition of $s_{m}$ yields $s_{m}=k_{1}+k_{2}+\cdots+k_{m}$.

Each $i\in\left\{  1,2,\ldots,m\right\}  $ satisfies%
\[
\dbinom{s_{i}}{k_{i}}=s_{i}!\cdot\dfrac{1}{k_{i}!}\cdot\dfrac{1}{s_{i-1}!}%
\]
\footnote{\textit{Proof:} Let $i\in\left\{  1,2,\ldots,m\right\}  $. Then,
both $i-1$ and $i$ belong to the set $\left\{  0,1,\ldots,m\right\}  $. Hence,
the definition of $s_{i-1}$ yields $s_{i-1}=k_{1}+k_{2}+\cdots+k_{i-1}%
\in\mathbb{N}$ (since $k_{1},k_{2},\ldots,k_{i-1}$ are elements of
$\mathbb{N}$). Thus, $s_{i-1}\geq0$. Meanwhile, the definition of $s_{i}$
yields
\[
s_{i}=k_{1}+k_{2}+\cdots+k_{i}=\underbrace{\left(  k_{1}+k_{2}+\cdots
+k_{i-1}\right)  }_{=s_{i-1}}+k_{i}=\underbrace{s_{i-1}}_{\geq0}+k_{i}\geq
k_{i}.
\]
Also, $s_{i}=k_{1}+k_{2}+\cdots+k_{i}\in\mathbb{N}$ (since $k_{1},k_{2}%
,\ldots,k_{i}$ are elements of $\mathbb{N}$). From $s_{i}=s_{i-1}+k_{i}$, we
obtain $s_{i}-k_{i}=s_{i-1}$.
\par
Now, (\ref{eq.binom.formula}) (applied to $s_{i}$ and $k_{i}$ instead of $m$
and $n$) yields
\begin{align*}
\dbinom{s_{i}}{k_{i}}  &  =\dfrac{s_{i}!}{k_{i}!\left(  s_{i}-k_{i}\right)
!}=\dfrac{s_{i}!}{k_{i}!s_{i-1}!}\ \ \ \ \ \ \ \ \ \ \left(  \text{since
}s_{i}-k_{i}=s_{i-1}\right) \\
&  =s_{i}!\cdot\dfrac{1}{k_{i}!}\cdot\dfrac{1}{s_{i-1}!}.
\end{align*}
}. Hence,%
\begin{align*}
\prod_{i=1}^{m}\underbrace{\dbinom{s_{i}}{k_{i}}}_{=s_{i}!\cdot\dfrac{1}%
{k_{i}!}\cdot\dfrac{1}{s_{i-1}!}}  &  =\prod_{i=1}^{m}\left(  s_{i}%
!\cdot\dfrac{1}{k_{i}!}\cdot\dfrac{1}{s_{i-1}!}\right)  =\underbrace{\left(
\prod_{i=1}^{m}s_{i}!\right)  }_{\substack{=\left(  \prod_{i=1}^{m-1}%
s_{i}!\right)  s_{m}!\\\text{(since }m\geq1\text{)}}}\cdot\left(  \prod
_{i=1}^{m}\dfrac{1}{k_{i}!}\right)  \cdot\underbrace{\left(  \prod_{i=1}%
^{m}\dfrac{1}{s_{i-1}!}\right)  }_{\substack{=\dfrac{1}{\prod_{i=1}^{m-1}%
s_{i}!}\\\text{(by (\ref{sol.multinom1.short.goal.pf.2}))}}}\\
&  =\left(  \prod_{i=1}^{m-1}s_{i}!\right)  s_{m}!\cdot\left(  \prod_{i=1}%
^{m}\dfrac{1}{k_{i}!}\right)  \cdot\dfrac{1}{\prod_{i=1}^{m-1}s_{i}!}\\
&  =\underbrace{s_{m}}_{=k_{1}+k_{2}+\cdots+k_{m}}!\cdot\underbrace{\left(
\prod_{i=1}^{m}\dfrac{1}{k_{i}!}\right)  }_{=\dfrac{1}{\prod_{i=1}^{m}k_{i}%
!}=\dfrac{1}{k_{1}!k_{2}!\cdots k_{m}!}}\\
&  =\left(  k_{1}+k_{2}+\cdots+k_{m}\right)  !\cdot\dfrac{1}{k_{1}%
!k_{2}!\cdots k_{m}!}=\dfrac{\left(  k_{1}+k_{2}+\cdots+k_{m}\right)  !}%
{k_{1}!k_{2}!\cdots k_{m}!}.
\end{align*}
This proves (\ref{sol.multinom1.short.goal}).]

On the other hand, each $i\in\left\{  1,2,\ldots,m\right\}  $ satisfies
$\dbinom{s_{i}}{k_{i}}\in\mathbb{Z}$ (by Proposition \ref{prop.binom.int}
(applied to $s_{i}$ and $k_{i}$ instead of $m$ and $n$)). In other words, for
each $i\in\left\{  1,2,\ldots,m\right\}  $, the number $\dbinom{s_{i}}{k_{i}}$
is an integer. Hence, $\prod_{i=1}^{m}\dbinom{s_{i}}{k_{i}}$ is an integer
(since a product of integers always is an integer). In other words,
$\dfrac{\left(  k_{1}+k_{2}+\cdots+k_{m}\right)  !}{k_{1}!k_{2}!\cdots k_{m}%
!}$ is an integer (since $\dfrac{\left(  k_{1}+k_{2}+\cdots+k_{m}\right)
!}{k_{1}!k_{2}!\cdots k_{m}!}=\prod_{i=1}^{m}\dbinom{s_{i}}{k_{i}}$). Since
$\dfrac{\left(  k_{1}+k_{2}+\cdots+k_{m}\right)  !}{k_{1}!k_{2}!\cdots k_{m}%
!}$ is clearly positive (because the numbers $\left(  k_{1}+k_{2}+\cdots
+k_{m}\right)  !$, $k_{1}!$, $k_{2}!$, $\ldots$, $k_{m}!$ are all positive),
we can thus conclude that $\dfrac{\left(  k_{1}+k_{2}+\cdots+k_{m}\right)
!}{k_{1}!k_{2}!\cdots k_{m}!}$ is a positive integer. This solves Exercise
\ref{exe.multinom1}.
\end{proof}
\end{vershort}

\begin{verlong}
\begin{proof}
[Solution to Exercise \ref{exe.multinom1}.]For each $i\in\left\{
0,1,\ldots,m\right\}  $, define an integer $s_{i}\in\mathbb{N}$ by
$s_{i}=k_{1}+k_{2}+\cdots+k_{i}$. We shall prove that%
\begin{equation}
\dfrac{\left(  k_{1}+k_{2}+\cdots+k_{m}\right)  !}{k_{1}!k_{2}!\cdots k_{m}%
!}=\prod_{i=1}^{m}\dbinom{s_{i}}{k_{i}}. \label{sol.multinom1.goal}%
\end{equation}

[\textit{Proof of (\ref{sol.multinom1.goal}):} If $m=0$, then
(\ref{sol.multinom1.goal}) holds\footnote{\textit{Proof.} Assume that $m=0$.
We must prove that (\ref{sol.multinom1.goal}) holds.
\par
Dividing the equality $\left(  \underbrace{k_{1}+k_{2}+\cdots+k_{0}}_{=\left(
\text{empty sum}\right)  =0}\right)  !=0!=1$ by the equality $k_{1}%
!k_{2}!\cdots k_{0}!=\left(  \text{empty product}\right)  =1$, we obtain
\[
\dfrac{\left(  k_{1}+k_{2}+\cdots+k_{0}\right)  !}{k_{1}!k_{2}!\cdots k_{0}%
!}=\dfrac{1}{1}=1.
\]
Comparing this with $\prod_{i=1}^{0}\dbinom{s_{i}}{k_{i}}=\left(  \text{empty
product}\right)  =1$, we find%
\[
\dfrac{\left(  k_{1}+k_{2}+\cdots+k_{0}\right)  !}{k_{1}!k_{2}!\cdots k_{0}%
!}=\prod_{i=1}^{0}\dbinom{s_{i}}{k_{i}}.
\]
This rewrites as $\dfrac{\left(  k_{1}+k_{2}+\cdots+k_{m}\right)  !}%
{k_{1}!k_{2}!\cdots k_{m}!}=\prod_{i=1}^{m}\dbinom{s_{i}}{k_{i}}$ (since
$m=0$). Thus, (\ref{sol.multinom1.goal}) holds. Qed.}. Hence, for the rest of
the proof of (\ref{sol.multinom1.goal}), we can WLOG assume that we don't have
$m=0$. Assume this.

We have $m\neq0$ (since we don't have $m=0$). Thus, $m\geq1$ (since
$m\in\mathbb{N}$), so that $m-1\geq0$.

We have $0\in\left\{  0,1,\ldots,m\right\}  $ (since $m\in\mathbb{N}$). The
definition of $s_{0}$ thus yields $s_{0}=k_{1}+k_{2}+\cdots+k_{0}=\left(
\text{empty sum}\right)  =0$. Hence, $s_{0}!=0!=1$. Now,%
\begin{align}
\prod_{i=1}^{m}\dfrac{1}{s_{i-1}!}  &  =\prod_{i=0}^{m-1}\dfrac{1}{s_{i}%
!}\ \ \ \ \ \ \ \ \ \ \left(  \text{here, we have substituted }i\text{ for
}i-1\text{ in the product}\right) \nonumber\\
&  =\underbrace{\dfrac{1}{s_{0}!}}_{\substack{=1\\\text{(since }%
s_{0}!=1\text{)}}}\prod_{i=1}^{m-1}\dfrac{1}{s_{i}!}%
\ \ \ \ \ \ \ \ \ \ \left(
\begin{array}
[c]{c}%
\text{here, we have split off the factor for }i=0\\
\text{from the product (since }m-1\geq0\text{)}%
\end{array}
\right) \nonumber\\
&  =\prod_{i=1}^{m-1}\dfrac{1}{s_{i}!}=\dfrac{1}{\prod_{i=1}^{m-1}s_{i}!}.
\label{sol.multinom1.goal.pf.2}%
\end{align}

We have $m\in\left\{  0,1,\ldots,m\right\}  $ (since $m\in\mathbb{N}$). The
definition of $s_{m}$ thus yields $s_{m}=k_{1}+k_{2}+\cdots+k_{m}$.

Each $i\in\left\{  1,2,\ldots,m\right\}  $ satisfies%
\[
\dbinom{s_{i}}{k_{i}}=s_{i}!\cdot\dfrac{1}{k_{i}!}\cdot\dfrac{1}{s_{i-1}!}%
\]
\footnote{\textit{Proof:} Let $i\in\left\{  1,2,\ldots,m\right\}  $. Then,
$i-1\in\left\{  0,1,\ldots,m-1\right\}  \subseteq\left\{  0,1,\ldots
,m\right\}  $. Hence, the definition of $s_{i-1}$ yields $s_{i-1}=k_{1}%
+k_{2}+\cdots+k_{i-1}\in\mathbb{N}$ (since $k_{1},k_{2},\ldots,k_{i-1}$ are
elements of $\mathbb{N}$). Thus, $s_{i-1}\geq0$. But $i\in\left\{
1,2,\ldots,m\right\}  \subseteq\left\{  0,1,\ldots,m\right\}  $. Hence, the
definition of $s_{i}$ yields
\[
s_{i}=k_{1}+k_{2}+\cdots+k_{i}=\underbrace{\left(  k_{1}+k_{2}+\cdots
+k_{i-1}\right)  }_{=s_{i-1}}+k_{i}=\underbrace{s_{i-1}}_{\geq0}+k_{i}\geq
k_{i}.
\]
Also, $s_{i}=k_{1}+k_{2}+\cdots+k_{i}\in\mathbb{N}$ (since $k_{1},k_{2}%
,\ldots,k_{i}$ are elements of $\mathbb{N}$). From $s_{i}=s_{i-1}+k_{i}$, we
obtain $s_{i}-k_{i}=s_{i-1}$.
\par
Now, (\ref{eq.binom.formula}) (applied to $s_{i}$ and $k_{i}$ instead of $m$
and $n$) yields
\begin{align*}
\dbinom{s_{i}}{k_{i}}  &  =\dfrac{s_{i}!}{k_{i}!\left(  s_{i}-k_{i}\right)
!}=\dfrac{s_{i}!}{k_{i}!s_{i-1}!}\ \ \ \ \ \ \ \ \ \ \left(  \text{since
}s_{i}-k_{i}=s_{i-1}\right) \\
&  =s_{i}!\cdot\dfrac{1}{k_{i}!}\cdot\dfrac{1}{s_{i-1}!}.
\end{align*}
Qed.}. Hence,%
\begin{align*}
&  \prod_{i=1}^{m}\underbrace{\dbinom{s_{i}}{k_{i}}}_{=s_{i}!\cdot\dfrac
{1}{k_{i}!}\cdot\dfrac{1}{s_{i-1}!}}\\
&  =\prod_{i=1}^{m}\left(  s_{i}!\cdot\dfrac{1}{k_{i}!}\cdot\dfrac{1}%
{s_{i-1}!}\right)  =\underbrace{\left(  \prod_{i=1}^{m}s_{i}!\right)
}_{\substack{=\left(  \prod_{i=1}^{m-1}s_{i}!\right)  s_{m}!\\\text{(here, we
have split off the}\\\text{factor for }i=m\text{ from the}\\\text{product
(since }m\geq1\text{))}}}\cdot\left(  \prod_{i=1}^{m}\dfrac{1}{k_{i}!}\right)
\cdot\underbrace{\left(  \prod_{i=1}^{m}\dfrac{1}{s_{i-1}!}\right)
}_{\substack{=\dfrac{1}{\prod_{i=1}^{m-1}s_{i}!}\\\text{(by
(\ref{sol.multinom1.goal.pf.2}))}}}\\
&  =\left(  \prod_{i=1}^{m-1}s_{i}!\right)  s_{m}!\cdot\left(  \prod_{i=1}%
^{m}\dfrac{1}{k_{i}!}\right)  \cdot\dfrac{1}{\prod_{i=1}^{m-1}s_{i}!}\\
&  =\underbrace{s_{m}}_{=k_{1}+k_{2}+\cdots+k_{m}}!\cdot\underbrace{\left(
\prod_{i=1}^{m}\dfrac{1}{k_{i}!}\right)  }_{=\dfrac{1}{\prod_{i=1}^{m}k_{i}%
!}=\dfrac{1}{k_{1}!k_{2}!\cdots k_{m}!}}\\
&  =\left(  k_{1}+k_{2}+\cdots+k_{m}\right)  !\cdot\dfrac{1}{k_{1}%
!k_{2}!\cdots k_{m}!}=\dfrac{\left(  k_{1}+k_{2}+\cdots+k_{m}\right)  !}%
{k_{1}!k_{2}!\cdots k_{m}!}.
\end{align*}
This proves (\ref{sol.multinom1.goal}).]

On the other hand, each $i\in\left\{  1,2,\ldots,m\right\}  $ satisfies
$\dbinom{s_{i}}{k_{i}}\in\mathbb{Z}$ (by Proposition \ref{prop.binom.int}
(applied to $s_{i}$ and $k_{i}$ instead of $m$ and $n$)). In other words, for
each $i\in\left\{  1,2,\ldots,m\right\}  $, the number $\dbinom{s_{i}}{k_{i}}$
is an integer. Hence, $\prod_{i=1}^{m}\dbinom{s_{i}}{k_{i}}$ is a product of
integers. Thus, $\prod_{i=1}^{m}\dbinom{s_{i}}{k_{i}}$ is an integer (since a
product of integers always is an integer). In other words, $\dfrac{\left(
k_{1}+k_{2}+\cdots+k_{m}\right)  !}{k_{1}!k_{2}!\cdots k_{m}!}$ is an integer
(since $\dfrac{\left(  k_{1}+k_{2}+\cdots+k_{m}\right)  !}{k_{1}!k_{2}!\cdots
k_{m}!}=\prod_{i=1}^{m}\dbinom{s_{i}}{k_{i}}$). Since $\dfrac{\left(
k_{1}+k_{2}+\cdots+k_{m}\right)  !}{k_{1}!k_{2}!\cdots k_{m}!}$ is clearly
positive (because the numbers $\left(  k_{1}+k_{2}+\cdots+k_{m}\right)  !$,
$k_{1}!$, $k_{2}!$, $\ldots$, $k_{m}!$ are all positive), we can thus conclude
that $\dfrac{\left(  k_{1}+k_{2}+\cdots+k_{m}\right)  !}{k_{1}!k_{2}!\cdots
k_{m}!}$ is a positive integer. This solves Exercise \ref{exe.multinom1}.
\end{proof}
\end{verlong}

\subsection{Solution to Exercise \ref{exe.bin.-1/2}}

\subsubsection{The solution}

We are going to give two proofs for each part of Exercise \ref{exe.bin.-1/2}:
one by direct manipulation of products, and one by induction on $n$. The
induction proofs will rely on the following fact:

\begin{lemma}
\label{lem.sol.bin.-1/2.m-1}Let $m$ be a positive integer. Let $q\in
\mathbb{Q}$. Then,%
\[
\dbinom{q}{m}=\dfrac{q-m+1}{m}\dbinom{q}{m-1}.
\]

\end{lemma}

\begin{proof}
[Proof of Lemma \ref{lem.sol.bin.-1/2.m-1}.]We have $m-1\in\mathbb{N}$ (since
$m$ is a positive integer). Thus, (\ref{eq.binom.mn}) (applied to $q$ and
$m-1$ instead of $m$ and $n$) yields%
\begin{align}
\dbinom{q}{m-1}  &  =\dfrac{q\left(  q-1\right)  \cdots\left(  q-\left(
m-1\right)  +1\right)  }{\left(  m-1\right)  !}\nonumber\\
&  =\dfrac{q\left(  q-1\right)  \cdots\left(  q-m+2\right)  }{\left(
m-1\right)  !} \label{pf.lem.sol.bin.-1/2.m-1.1}%
\end{align}
(since $q-\left(  m-1\right)  +1=q-m+2$).

On the other hand, (\ref{eq.n!.rec}) (applied to $n=m$) yields $m!=m\cdot
\left(  m-1\right)  !$. Also, (\ref{eq.binom.mn}) (applied to $q$ and $m$
instead of $m$ and $n$) yields%
\begin{align*}
\dbinom{q}{m}  &  =\dfrac{q\left(  q-1\right)  \cdots\left(  q-m+1\right)
}{m!}=\underbrace{\dfrac{1}{m!}}_{\substack{=\dfrac{1}{m\cdot\left(
m-1\right)  !}\\\text{(since }m!=m\cdot\left(  m-1\right)  !\text{)}%
}}\underbrace{\left(  q\left(  q-1\right)  \cdots\left(  q-m+1\right)
\right)  }_{\substack{=\left(  q\left(  q-1\right)  \cdots\left(
q-m+2\right)  \right)  \cdot\left(  q-m+1\right)  \\\text{(since }m\text{ is a
positive integer)}}}\\
&  =\dfrac{1}{m\cdot\left(  m-1\right)  !}\left(  q\left(  q-1\right)
\cdots\left(  q-m+2\right)  \right)  \cdot\left(  q-m+1\right) \\
&  =\dfrac{q-m+1}{m}\cdot\underbrace{\dfrac{q\left(  q-1\right)  \cdots\left(
q-m+2\right)  }{\left(  m-1\right)  !}}_{\substack{=\dbinom{q}{m-1}\\\text{(by
(\ref{pf.lem.sol.bin.-1/2.m-1.1}))}}}=\dfrac{q-m+1}{m}\dbinom{q}{m-1}.
\end{align*}
This proves Lemma \ref{lem.sol.bin.-1/2.m-1}.
\end{proof}

We shall also use the following simple observation in the solutions to
Exercise \ref{exe.bin.-1/2} \textbf{(b)}:

\begin{lemma}
\label{lem.sol.bin.-1/2.2nn}Let $n\in\mathbb{N}$. Then,%
\[
\dbinom{2n}{n}=\dfrac{\left(  2n\right)  !}{n!^{2}}.
\]

\end{lemma}

\begin{proof}
[Proof of Lemma \ref{lem.sol.bin.-1/2.2nn}.]From $n\in\mathbb{N}$, we obtain
$n\geq0$, and thus $2n\geq n\geq0$. Hence, $2n\in\mathbb{N}$. Thus,
Proposition \ref{prop.binom.formula} (applied to $m=2n$) yields%
\[
\dbinom{2n}{n}=\dfrac{\left(  2n\right)  !}{n!\left(  2n-n\right)  !}%
=\dfrac{\left(  2n\right)  !}{n!n!}%
\]
(since $2n-n=n$). Thus, $\dbinom{2n}{n}=\dfrac{\left(  2n\right)  !}%
{n!n!}=\dfrac{\left(  2n\right)  !}{n!^{2}}$. This proves Lemma
\ref{lem.sol.bin.-1/2.2nn}.
\end{proof}

\begin{proof}
[Solution to Exercise \ref{exe.bin.-1/2}.]\textbf{(a)} \textit{First solution
to Exercise \ref{exe.bin.-1/2} \textbf{(a)}:} We have%
\begin{align*}
\left(  2n\right)  !  &  =1\cdot2\cdot\cdots\cdot\left(  2n\right)
=\prod_{k\in\left\{  1,2,\ldots,2n\right\}  }k=\underbrace{\left(
\prod_{\substack{k\in\left\{  1,2,\ldots,2n\right\}  ;\\k\text{ is even}%
}}k\right)  }_{\substack{=2\cdot4\cdot6\cdot\cdots\cdot\left(  2n\right)
\\=\prod_{i=1}^{n}\left(  2i\right)  \\=2^{n}\prod_{i=1}^{n}i}%
}\underbrace{\left(  \prod_{\substack{k\in\left\{  1,2,\ldots,2n\right\}
;\\k\text{ is odd}}}k\right)  }_{=1\cdot3\cdot5\cdot\cdots\cdot\left(
2n-1\right)  }\\
&  \ \ \ \ \ \ \ \ \ \ \left(
\begin{array}
[c]{c}%
\text{here, we have split the product }\prod_{k\in\left\{  1,2,\ldots
,2n\right\}  }k\text{ into one product}\\
\text{containing all even }k\text{ and one product containing all odd }k
\end{array}
\right) \\
&  =2^{n}\underbrace{\left(  \prod_{i=1}^{n}i\right)  }_{=1\cdot2\cdot
\cdots\cdot n=n!}\cdot\left(  1\cdot3\cdot5\cdot\cdots\cdot\left(
2n-1\right)  \right) \\
&  =2^{n}n!\cdot\left(  1\cdot3\cdot5\cdot\cdots\cdot\left(  2n-1\right)
\right)  .
\end{align*}
Dividing this equality by $2^{n}n!$, we obtain%
\[
\dfrac{\left(  2n\right)  !}{2^{n}n!}=1\cdot3\cdot5\cdot\cdots\cdot\left(
2n-1\right)  =\left(  2n-1\right)  \cdot\left(  2n-3\right)  \cdot\cdots
\cdot1.
\]
This solves Exercise \ref{exe.bin.-1/2} \textbf{(a)}.

\textit{Second solution to Exercise \ref{exe.bin.-1/2} \textbf{(a)}:} Let us
solve Exercise \ref{exe.bin.-1/2} \textbf{(a)} by induction on $n$:

\textit{Induction base:} From $2\cdot0=0$, we obtain%
\begin{align*}
\dfrac{\left(  2\cdot0\right)  !}{2^{0}0!}  &  =\dfrac{0!}{2^{0}0!}=\dfrac
{1}{2^{0}}=\dfrac{1}{1}\ \ \ \ \ \ \ \ \ \ \left(  \text{since }2^{0}=1\right)
\\
&  =1.
\end{align*}
Comparing this with%
\[
\left(  2\cdot0-1\right)  \cdot\left(  2\cdot0-3\right)  \cdot\cdots
\cdot1=\left(  \text{empty product}\right)  =1,
\]
we find $\left(  2\cdot0-1\right)  \cdot\left(  2\cdot0-3\right)  \cdot
\cdots\cdot1=\dfrac{\left(  2\cdot0\right)  !}{2^{0}0!}$. In other words,
Exercise \ref{exe.bin.-1/2} \textbf{(a)} holds for $n=0$. This completes the
induction base.

\textit{Induction step:} Let $m\in\mathbb{N}$. Assume that Exercise
\ref{exe.bin.-1/2} \textbf{(a)} holds for $n=m$. We must prove that Exercise
\ref{exe.bin.-1/2} \textbf{(a)} holds for $n=m+1$.

We have assumed that Exercise \ref{exe.bin.-1/2} \textbf{(a)} holds for $n=m$.
In other words, we have%
\begin{equation}
\left(  2m-1\right)  \cdot\left(  2m-3\right)  \cdot\cdots\cdot1=\dfrac
{\left(  2m\right)  !}{2^{m}m!}. \label{sol.bin.-1/2.a.IH}%
\end{equation}

Clearly, $m+1$ is a positive integer (since $m\in\mathbb{N}$). Hence,
(\ref{eq.n!.rec}) (applied to $n=m+1$) yields $\left(  m+1\right)  !=\left(
m+1\right)  \cdot\left(  \underbrace{\left(  m+1\right)  -1}_{=m}\right)
!=\left(  m+1\right)  \cdot m!$.

But $2\left(  m+1\right)  =2m+2$, and thus%
\begin{align*}
\left(  2\left(  m+1\right)  \right)  !  &  =\left(  2m+2\right)
!=1\cdot2\cdot\cdots\cdot\left(  2m+2\right) \\
&  =\underbrace{\left(  1\cdot2\cdot\cdots\cdot\left(  2m\right)  \right)
}_{\substack{=\left(  2m\right)  !\\\text{(since }\left(  2m\right)
!=1\cdot2\cdot\cdots\cdot\left(  2m\right)  \text{)}}}\cdot\left(
2m+1\right)  \cdot\left(  2m+2\right) \\
&  =\left(  2m\right)  !\cdot\left(  2m+1\right)  \cdot\left(  2m+2\right)  .
\end{align*}
Now,%
\begin{align*}
\dfrac{\left(  2\left(  m+1\right)  \right)  !}{2^{m+1}\left(  m+1\right)  !}
&  =\underbrace{\left(  2\left(  m+1\right)  \right)  !}_{=\left(  2m\right)
!\cdot\left(  2m+1\right)  \cdot\left(  2m+2\right)  }/\left(
\underbrace{2^{m+1}}_{=2\cdot2^{m}}\underbrace{\left(  m+1\right)
!}_{=\left(  m+1\right)  \cdot m!}\right) \\
&  =\left(  2m\right)  !\cdot\left(  2m+1\right)  \cdot\left(  2m+2\right)
/\left(  2\cdot2^{m}\cdot\left(  m+1\right)  \cdot m!\right) \\
&  =\dfrac{\left(  2m\right)  !\cdot\left(  2m+1\right)  \cdot\left(
2m+2\right)  }{2\cdot2^{m}\cdot\left(  m+1\right)  \cdot m!}=\dfrac{\left(
2m\right)  !}{2^{m}m!}\cdot\left(  2m+1\right)  \cdot\underbrace{\dfrac
{2m+2}{2\left(  m+1\right)  }}_{=1}\\
&  =\dfrac{\left(  2m\right)  !}{2^{m}m!}\cdot\left(  2m+1\right)  =\left(
2m+1\right)  \cdot\dfrac{\left(  2m\right)  !}{2^{m}m!}.
\end{align*}
Comparing this with%
\begin{align*}
&  \left(  2\left(  m+1\right)  -1\right)  \cdot\left(  2\left(  m+1\right)
-3\right)  \cdot\cdots\cdot1\\
&  =\underbrace{\left(  2\left(  m+1\right)  -1\right)  }_{=2m+1}%
\cdot\underbrace{\left(  \left(  2\left(  m+1\right)  -3\right)  \cdot\left(
2\left(  m+1\right)  -5\right)  \cdot\cdots\cdot1\right)  }%
_{\substack{=\left(  2m-1\right)  \cdot\left(  2m-3\right)  \cdot\cdots
\cdot1\\\text{(since }2\left(  m+1\right)  -3=2m-1\text{)}}}\\
&  =\left(  2m+1\right)  \cdot\underbrace{\left(  \left(  2m-1\right)
\cdot\left(  2m-3\right)  \cdot\cdots\cdot1\right)  }_{\substack{=\dfrac
{\left(  2m\right)  !}{2^{m}m!}\\\text{(by (\ref{sol.bin.-1/2.a.IH}))}%
}}=\left(  2m+1\right)  \cdot\dfrac{\left(  2m\right)  !}{2^{m}m!},
\end{align*}
we obtain%
\[
\left(  2\left(  m+1\right)  -1\right)  \cdot\left(  2\left(  m+1\right)
-3\right)  \cdot\cdots\cdot1=\dfrac{\left(  2\left(  m+1\right)  \right)
!}{2^{m+1}\left(  m+1\right)  !}.
\]
In other words, Exercise \ref{exe.bin.-1/2} \textbf{(a)} holds for $n=m+1$.
This completes the induction step. Thus, Exercise \ref{exe.bin.-1/2}
\textbf{(a)} is solved.

\textbf{(b)} \textit{First solution to Exercise \ref{exe.bin.-1/2}
\textbf{(b)}:} The equality (\ref{eq.binom.mn}) (applied to $-1/2$ instead of
$m$) yields%
\begin{align*}
\dbinom{-1/2}{n}  &  =\dfrac{\left(  -1/2\right)  \left(  -1/2-1\right)
\cdots\left(  -1/2-n+1\right)  }{n!}\\
&  =\dfrac{1}{n!}\cdot\underbrace{\left(  -1/2\right)  \left(  -1/2-1\right)
\cdots\left(  -1/2-n+1\right)  }_{=\prod_{k=0}^{n-1}\left(  -1/2-k\right)  }\\
&  =\dfrac{1}{n!}\cdot\prod_{k=0}^{n-1}\underbrace{\left(  -1/2-k\right)
}_{=\dfrac{-1}{2}\left(  2k+1\right)  }=\dfrac{1}{n!}\cdot\underbrace{\prod
_{k=0}^{n-1}\left(  \dfrac{-1}{2}\left(  2k+1\right)  \right)  }_{=\left(
\dfrac{-1}{2}\right)  ^{n}\prod_{k=0}^{n-1}\left(  2k+1\right)  }\\
&  =\dfrac{1}{n!}\cdot\left(  \dfrac{-1}{2}\right)  ^{n}\underbrace{\prod
_{k=0}^{n-1}\left(  2k+1\right)  }_{\substack{=1\cdot3\cdot5\cdot\cdots
\cdot\left(  2\left(  n-1\right)  +1\right)  \\=1\cdot3\cdot5\cdot\cdots
\cdot\left(  2n-1\right)  \\\text{(since }2\left(  n-1\right)  +1=2n-1\text{)}%
}}=\dfrac{1}{n!}\cdot\left(  \dfrac{-1}{2}\right)  ^{n}\underbrace{\left(
1\cdot3\cdot5\cdot\cdots\cdot\left(  2n-1\right)  \right)  }%
_{\substack{=\left(  2n-1\right)  \cdot\left(  2n-3\right)  \cdot\cdots
\cdot1\\=\dfrac{\left(  2n\right)  !}{2^{n}n!}\\\text{(by Exercise
\ref{exe.bin.-1/2} \textbf{(a)})}}}\\
&  =\dfrac{1}{n!}\cdot\left(  \dfrac{-1}{2}\right)  ^{n}\cdot\dfrac{\left(
2n\right)  !}{2^{n}n!}=\left(  \dfrac{-1}{2}\right)  ^{n}\cdot
\underbrace{\dfrac{1}{2^{n}}}_{=\left(  \dfrac{1}{2}\right)  ^{n}}\cdot
\dfrac{\left(  2n\right)  !}{n!^{2}}=\left(  \dfrac{-1}{2}\right)  ^{n}%
\cdot\left(  \dfrac{1}{2}\right)  ^{n}\cdot\dfrac{\left(  2n\right)  !}%
{n!^{2}}.
\end{align*}
Comparing this with%
\[
\left(  \underbrace{\dfrac{-1}{4}}_{=\dfrac{-1}{2}\cdot\dfrac{1}{2}}\right)
^{n}\underbrace{\dbinom{2n}{n}}_{\substack{=\dfrac{\left(  2n\right)
!}{n!^{2}}\\\text{(by Lemma \ref{lem.sol.bin.-1/2.2nn})}}}=\underbrace{\left(
\dfrac{-1}{2}\cdot\dfrac{1}{2}\right)  ^{n}}_{=\left(  \dfrac{-1}{2}\right)
^{n}\cdot\left(  \dfrac{1}{2}\right)  ^{n}}\cdot\dfrac{\left(  2n\right)
!}{n!^{2}}=\left(  \dfrac{-1}{2}\right)  ^{n}\cdot\left(  \dfrac{1}{2}\right)
^{n}\cdot\dfrac{\left(  2n\right)  !}{n!^{2}},
\]
we obtain $\dbinom{-1/2}{n}=\left(  \dfrac{-1}{4}\right)  ^{n}\dbinom{2n}{n}$.
This solves Exercise \ref{exe.bin.-1/2} \textbf{(b)}.

\textit{Second solution to Exercise \ref{exe.bin.-1/2} \textbf{(b)}:} Let us
solve Exercise \ref{exe.bin.-1/2} \textbf{(b)} by induction on $n$:

\textit{Induction base:} Applying (\ref{eq.binom.00}) to $m=2\cdot0$, we
obtain $\dbinom{2\cdot0}{0}=1$. But Proposition \ref{prop.binom.00}
\textbf{(a)} (applied to $-1/2$ instead of $m$) yields $\dbinom{-1/2}{0}=1$.
Comparing this with $\underbrace{\left(  \dfrac{-1}{4}\right)  ^{0}}%
_{=1}\underbrace{\dbinom{2\cdot0}{0}}_{=1}=1$, we obtain $\dbinom{-1/2}%
{0}=\left(  \dfrac{-1}{4}\right)  ^{0}\dbinom{2\cdot0}{0}$. In other words,
Exercise \ref{exe.bin.-1/2} \textbf{(b)} holds for $n=0$. This completes the
induction base.

\textit{Induction step:} Let $m$ be a positive integer. Assume that Exercise
\ref{exe.bin.-1/2} \textbf{(b)} holds for $n=m-1$. We must now prove that
Exercise \ref{exe.bin.-1/2} \textbf{(b)} holds for $n=m$.

We have assumed that Exercise \ref{exe.bin.-1/2} \textbf{(b)} holds for
$n=m-1$. In other words, we have%
\begin{equation}
\dbinom{-1/2}{m-1}=\left(  \dfrac{-1}{4}\right)  ^{m-1}\dbinom{2\left(
m-1\right)  }{m-1}. \label{sol.bin.-1/2.2nd.IH}%
\end{equation}

We have $m-1\in\mathbb{N}$ (since $m$ is a positive integer). Thus, Lemma
\ref{lem.sol.bin.-1/2.2nn} (applied to $n=m-1$) yields
\begin{equation}
\dbinom{2\left(  m-1\right)  }{m-1}=\dfrac{\left(  2\left(  m-1\right)
\right)  !}{\left(  m-1\right)  !^{2}}. \label{sol.bin.-1/2.2nd.3}%
\end{equation}

Applying (\ref{eq.n!.rec}) to $n=m$, we obtain $m!=m\cdot\left(  m-1\right)
!$. Also, $m\geq1$ (since $m$ is a positive integer), so that $2m\geq2$ and
thus $2m-1\geq2-1=1$. Hence, $2m-1$ is a positive integer. Thus,
(\ref{eq.n!.rec}) (applied to $n=2m-1$) yields
\[
\left(  2m-1\right)  !=\left(  2m-1\right)  \cdot\left(  \underbrace{\left(
2m-1\right)  -1}_{=2\left(  m-1\right)  }\right)  !=\left(  2m-1\right)
\cdot\left(  2\left(  m-1\right)  \right)  !.
\]
Moreover, $2m$ is a positive integer (since $2m\geq2>0$). Thus,
(\ref{eq.n!.rec}) (applied to $n=2m$) yields
\[
\left(  2m\right)  !=\left(  2m\right)  \cdot\underbrace{\left(  2m-1\right)
!}_{=\left(  2m-1\right)  \cdot\left(  2\left(  m-1\right)  \right)
!}=\left(  2m\right)  \cdot\left(  2m-1\right)  \cdot\left(  2\left(
m-1\right)  \right)  !.
\]
But Lemma \ref{lem.sol.bin.-1/2.2nn} (applied to $n=m$) yields%
\begin{align*}
\dbinom{2m}{m}  &  =\dfrac{\left(  2m\right)  !}{m!^{2}}=\underbrace{\left(
2m\right)  !}_{=\left(  2m\right)  \cdot\left(  2m-1\right)  \cdot\left(
2\left(  m-1\right)  \right)  !}/\left(  \underbrace{m!}_{=m\cdot\left(
m-1\right)  !}\right)  ^{2}\\
&  =\left(  2m\right)  \cdot\left(  2m-1\right)  \cdot\left(  2\left(
m-1\right)  \right)  !/\left(  m\cdot\left(  m-1\right)  !\right)  ^{2}\\
&  =\dfrac{\left(  2m\right)  \cdot\left(  2m-1\right)  \cdot\left(  2\left(
m-1\right)  \right)  !}{\left(  m\cdot\left(  m-1\right)  !\right)  ^{2}%
}=\underbrace{\dfrac{2m}{m}}_{=2}\cdot\dfrac{2m-1}{m}\cdot\underbrace{\dfrac
{\left(  2\left(  m-1\right)  \right)  !}{\left(  m-1\right)  !^{2}}%
}_{\substack{=\dbinom{2\left(  m-1\right)  }{m-1}\\\text{(by
(\ref{sol.bin.-1/2.2nd.3}))}}}\\
&  =2\cdot\dfrac{2m-1}{m}\cdot\dbinom{2\left(  m-1\right)  }{m-1}.
\end{align*}
Hence,%
\begin{align*}
&  \underbrace{\left(  \dfrac{-1}{4}\right)  ^{m}}_{=\left(  \dfrac{-1}%
{4}\right)  ^{m-1}\cdot\dfrac{-1}{4}}\ \ \ \underbrace{\dbinom{2m}{m}%
}_{=2\cdot\dfrac{2m-1}{m}\cdot\dbinom{2\left(  m-1\right)  }{m-1}}\\
=  &  \left(  \dfrac{-1}{4}\right)  ^{m-1}\cdot\underbrace{\dfrac{-1}{4}%
\cdot2}_{=\dfrac{-1}{2}}\cdot\dfrac{2m-1}{m}\cdot\dbinom{2\left(  m-1\right)
}{m-1}\\
=  &  \left(  \dfrac{-1}{4}\right)  ^{m-1}\cdot\dfrac{-1}{2}\cdot\dfrac
{2m-1}{m}\cdot\dbinom{2\left(  m-1\right)  }{m-1}.
\end{align*}
Comparing this with%
\begin{align*}
\dbinom{-1/2}{m}  &  =\underbrace{\dfrac{-1/2-m+1}{m}}_{=\dfrac{-1}{2}%
\cdot\dfrac{2m-1}{m}}\underbrace{\dbinom{-1/2}{m-1}}_{\substack{=\left(
\dfrac{-1}{4}\right)  ^{m-1}\dbinom{2\left(  m-1\right)  }{m-1}\\\text{(by
(\ref{sol.bin.-1/2.2nd.IH}))}}}\\
&  \ \ \ \ \ \ \ \ \ \ \left(  \text{by Lemma \ref{lem.sol.bin.-1/2.m-1}
(applied to }q=-1/2\text{)}\right) \\
&  =\dfrac{-1}{2}\cdot\dfrac{2m-1}{m}\cdot\left(  \dfrac{-1}{4}\right)
^{m-1}\dbinom{2\left(  m-1\right)  }{m-1}\\
&  =\left(  \dfrac{-1}{4}\right)  ^{m-1}\cdot\dfrac{-1}{2}\cdot\dfrac{2m-1}%
{m}\cdot\dbinom{2\left(  m-1\right)  }{m-1},
\end{align*}
we obtain $\dbinom{-1/2}{m}=\left(  \dfrac{-1}{4}\right)  ^{m}\dbinom{2m}{m}$.
In other words, Exercise \ref{exe.bin.-1/2} \textbf{(b)} holds for $n=m$. This
completes the induction step. Thus, the induction proof of Exercise
\ref{exe.bin.-1/2} \textbf{(b)} is complete.

\textbf{(c)} \textit{First solution to Exercise \ref{exe.bin.-1/2}
\textbf{(c)}:} We have
\begin{align*}
\left(  3n\right)  !  &  =1\cdot2\cdot\cdots\cdot\left(  3n\right)
=\prod_{k\in\left\{  1,2,\ldots,3n\right\}  }k\\
&  =\left(  \prod_{\substack{k\in\left\{  1,2,\ldots,3n\right\}
;\\k\equiv0\operatorname{mod}3}}k\right)  \left(  \prod_{\substack{k\in
\left\{  1,2,\ldots,3n\right\}  ;\\k\equiv1\operatorname{mod}3}}k\right)
\left(  \prod_{\substack{k\in\left\{  1,2,\ldots,3n\right\}  ;\\k\equiv
2\operatorname{mod}3}}k\right)
\end{align*}
(here, we have split the product $\prod_{k\in\left\{  1,2,\ldots,3n\right\}
}k$ into three smaller products, because each $k\in\left\{  1,2,\ldots
,3n\right\}  $ must satisfy exactly one of the three conditions $k\equiv
0\operatorname{mod}3$, $k\equiv1\operatorname{mod}3$ and $k\equiv
2\operatorname{mod}3$). Thus,%
\begin{align}
\left(  3n\right)  !  &  =\underbrace{\left(  \prod_{\substack{k\in\left\{
1,2,\ldots,3n\right\}  ;\\k\equiv0\operatorname{mod}3}}k\right)
}_{\substack{=3\cdot6\cdot9\cdot\cdots\cdot\left(  3n\right)  \\=\prod
_{i=1}^{n}\left(  3i\right)  \\=3^{n}\prod_{i=1}^{n}i}}\underbrace{\left(
\prod_{\substack{k\in\left\{  1,2,\ldots,3n\right\}  ;\\k\equiv
1\operatorname{mod}3}}k\right)  }_{\substack{=1\cdot4\cdot7\cdot\cdots
\cdot\left(  3n-2\right)  \\=\prod_{i=0}^{n-1}\left(  3i+1\right)
}}\underbrace{\left(  \prod_{\substack{k\in\left\{  1,2,\ldots,3n\right\}
;\\k\equiv2\operatorname{mod}3}}k\right)  }_{\substack{=2\cdot5\cdot
8\cdot\cdots\cdot\left(  3n-1\right)  \\=\prod_{i=0}^{n-1}\left(  3i+2\right)
}}\nonumber\\
&  =3^{n}\underbrace{\left(  \prod_{i=1}^{n}i\right)  }_{=n!}\left(
\prod_{i=0}^{n-1}\left(  3i+1\right)  \right)  \left(  \prod_{i=0}%
^{n-1}\left(  3i+2\right)  \right) \nonumber\\
&  =3^{n}n!\left(  \prod_{i=0}^{n-1}\left(  3i+1\right)  \right)  \left(
\prod_{i=0}^{n-1}\left(  3i+2\right)  \right)  . \label{sol.bin.-1/2.c.1}%
\end{align}

On the other hand, for each $g\in\mathbb{Z}$, we have%
\begin{align}
\dbinom{-g/3}{n}  &  =\dfrac{\left(  -g/3\right)  \left(  -g/3-1\right)
\cdots\left(  -g/3-n+1\right)  }{n!}\nonumber\\
&  \ \ \ \ \ \ \ \ \ \ \left(  \text{by (\ref{eq.binom.mn}) (applied to
}m=-g/3\text{)}\right) \nonumber\\
&  =\dfrac{1}{n!}\cdot\underbrace{\left(  -g/3\right)  \left(  -g/3-1\right)
\cdots\left(  -g/3-n+1\right)  }_{=\prod_{i=0}^{n-1}\left(  -g/3-i\right)
}\nonumber\\
&  =\dfrac{1}{n!}\cdot\prod_{i=0}^{n-1}\underbrace{\left(  -g/3-i\right)
}_{=\dfrac{-1}{3}\left(  3i+g\right)  }=\dfrac{1}{n!}\cdot\underbrace{\prod
_{i=0}^{n-1}\left(  \dfrac{-1}{3}\left(  3i+g\right)  \right)  }_{=\left(
\dfrac{-1}{3}\right)  ^{n}\prod_{i=0}^{n-1}\left(  3i+g\right)  }\nonumber\\
&  =\dfrac{1}{n!}\cdot\left(  \dfrac{-1}{3}\right)  ^{n}\prod_{i=0}%
^{n-1}\left(  3i+g\right)  . \label{sol.bin.-1/2.c.2}%
\end{align}

Now,%
\begin{align*}
&  \underbrace{\left(  3^{n}n!\right)  ^{3}}_{=3^{n}n!\cdot3^{n}n!\cdot
3^{n}n!}\ \ \underbrace{\dbinom{-1/3}{n}}_{\substack{=\dfrac{1}{n!}%
\cdot\left(  \dfrac{-1}{3}\right)  ^{n}\prod_{i=0}^{n-1}\left(  3i+1\right)
\\\text{(by (\ref{sol.bin.-1/2.c.2}), applied to }g=1\text{)}}%
}\ \ \underbrace{\dbinom{-2/3}{n}}_{\substack{=\dfrac{1}{n!}\cdot\left(
\dfrac{-1}{3}\right)  ^{n}\prod_{i=0}^{n-1}\left(  3i+2\right)  \\\text{(by
(\ref{sol.bin.-1/2.c.2}), applied to }g=2\text{)}}}\\
&  =3^{n}n!\cdot3^{n}n!\cdot3^{n}n!\cdot\left(  \dfrac{1}{n!}\cdot\left(
\dfrac{-1}{3}\right)  ^{n}\prod_{i=0}^{n-1}\left(  3i+1\right)  \right)
\cdot\left(  \dfrac{1}{n!}\cdot\left(  \dfrac{-1}{3}\right)  ^{n}\prod
_{i=0}^{n-1}\left(  3i+2\right)  \right) \\
&  =\underbrace{3^{n}\cdot3^{n}\cdot3^{n}\cdot\left(  \dfrac{-1}{3}\right)
^{n}\left(  \dfrac{-1}{3}\right)  ^{n}}_{=\left(  3\cdot3\cdot3\cdot\dfrac
{-1}{3}\cdot\dfrac{-1}{3}\right)  ^{n}}n!\left(  \prod_{i=0}^{n-1}\left(
3i+1\right)  \right)  \left(  \prod_{i=0}^{n-1}\left(  3i+2\right)  \right) \\
&  =\left(  \underbrace{3\cdot3\cdot3\cdot\dfrac{-1}{3}\cdot\dfrac{-1}{3}%
}_{=3}\right)  ^{n}n!\left(  \prod_{i=0}^{n-1}\left(  3i+1\right)  \right)
\left(  \prod_{i=0}^{n-1}\left(  3i+2\right)  \right) \\
&  =3^{n}n!\left(  \prod_{i=0}^{n-1}\left(  3i+1\right)  \right)  \left(
\prod_{i=0}^{n-1}\left(  3i+2\right)  \right)  .
\end{align*}
Comparing this with (\ref{sol.bin.-1/2.c.1}), we obtain $\left(
3^{n}n!\right)  ^{3}\dbinom{-1/3}{n}\dbinom{-2/3}{n}=\left(  3n\right)  !$.
Dividing this equality by $\left(  3^{n}n!\right)  ^{3}$, we find
$\dbinom{-1/3}{n}\dbinom{-2/3}{n}=\dfrac{\left(  3n\right)  !}{\left(
3^{n}n!\right)  ^{3}}$. This solves Exercise \ref{exe.bin.-1/2} \textbf{(c)}.

\textit{Second solution to Exercise \ref{exe.bin.-1/2} \textbf{(c)}:} Let us
solve Exercise \ref{exe.bin.-1/2} \textbf{(c)} by induction on $n$:

\textit{Induction base:} Proposition \ref{prop.binom.00} \textbf{(a)} (applied
to $m=-1/3$) yields $\dbinom{-1/3}{0}=1$. Also, Proposition
\ref{prop.binom.00} \textbf{(a)} (applied to $m=-2/3$) yields $\dbinom
{-2/3}{0}=1$. On the other hand,%
\[
\dfrac{\left(  3\cdot0\right)  !}{\left(  3^{0}0!\right)  ^{3}}=\left(
\underbrace{3\cdot0}_{=0}\right)  !/\left(  \underbrace{3^{0}}_{=1}%
\underbrace{0!}_{=1}\right)  ^{3}=0!/1^{3}=0!=1.
\]
Comparing this with%
\[
\underbrace{\dbinom{-1/3}{0}}_{=1}\underbrace{\dbinom{-2/3}{0}}_{=1}=1,
\]
we obtain $\dbinom{-1/3}{0}\dbinom{-2/3}{0}=\dfrac{\left(  3\cdot0\right)
!}{\left(  3^{0}0!\right)  ^{3}}$. In other words, Exercise \ref{exe.bin.-1/2}
\textbf{(c)} holds for $n=0$. This completes the induction base.

\textit{Induction step:} Let $m$ be a positive integer. Assume that Exercise
\ref{exe.bin.-1/2} \textbf{(c)} holds for $n=m-1$. We must now prove that
Exercise \ref{exe.bin.-1/2} \textbf{(c)} holds for $n=m$.

We have assumed that Exercise \ref{exe.bin.-1/2} \textbf{(c)} holds for
$n=m-1$. In other words, we have%
\begin{equation}
\dbinom{-1/3}{m-1}\dbinom{-2/3}{m-1}=\dfrac{\left(  3\left(  m-1\right)
\right)  !}{\left(  3^{m-1}\left(  m-1\right)  !\right)  ^{3}}.
\label{sol.bin.-1/2.c.2nd.IH}%
\end{equation}

We have $m-1\in\mathbb{N}$ (since $m$ is a positive integer). Hence, $3\left(
m-1\right)  \in\mathbb{N}$. The definition of $\left(  3\left(  m-1\right)
\right)  !$ thus yields%
\begin{equation}
\left(  3\left(  m-1\right)  \right)  !=1\cdot2\cdot\cdots\cdot\left(
3\left(  m-1\right)  \right)  =1\cdot2\cdot\cdots\cdot\left(  3m-3\right)
\label{sol.bin.-1/2.c.2nd.1}%
\end{equation}
(since $3\left(  m-1\right)  =3m-3$).

The definition of $\left(  3m\right)  !$ yields
\begin{align*}
\left(  3m\right)  !  &  =1\cdot2\cdot\cdots\cdot\left(  3m\right)
=\underbrace{\left(  1\cdot2\cdot\cdots\cdot\left(  3m-3\right)  \right)
}_{\substack{=\left(  3\left(  m-1\right)  \right)  !\\\text{(by
(\ref{sol.bin.-1/2.c.2nd.1}))}}}\cdot\left(  \left(  3m-2\right)  \cdot\left(
3m-1\right)  \cdot\left(  3m\right)  \right) \\
&  =\left(  3\left(  m-1\right)  \right)  !\cdot\left(  \left(  3m-2\right)
\cdot\left(  3m-1\right)  \cdot\left(  3m\right)  \right)  .
\end{align*}

Applying (\ref{eq.n!.rec}) to $n=m$, we obtain $m!=m\cdot\left(  m-1\right)
!$.

Now,%
\begin{align}
\dfrac{\left(  3m\right)  !}{\left(  3^{m}m!\right)  ^{3}}  &
=\underbrace{\left(  3m\right)  !}_{=\left(  3\left(  m-1\right)  \right)
!\cdot\left(  \left(  3m-2\right)  \cdot\left(  3m-1\right)  \cdot\left(
3m\right)  \right)  }/\left(  \underbrace{3^{m}}_{=3\cdot3^{m-1}%
}\underbrace{m!}_{=m\cdot\left(  m-1\right)  !}\right)  ^{3}\nonumber\\
&  =\left(  3\left(  m-1\right)  \right)  !\cdot\left(  \left(  3m-2\right)
\cdot\left(  3m-1\right)  \cdot\left(  3m\right)  \right)
/\underbrace{\left(  3\cdot3^{m-1}m\cdot\left(  m-1\right)  !\right)  ^{3}%
}_{=3^{3}m^{3}\left(  3^{m-1}\left(  m-1\right)  !\right)  ^{3}}\nonumber\\
&  =\left(  3\left(  m-1\right)  \right)  !\cdot\left(  \left(  3m-2\right)
\cdot\left(  3m-1\right)  \cdot\left(  3m\right)  \right)  /\left(  3^{3}%
m^{3}\left(  3^{m-1}\left(  m-1\right)  !\right)  ^{3}\right) \nonumber\\
&  =\dfrac{\left(  3\left(  m-1\right)  \right)  !\cdot\left(  \left(
3m-2\right)  \cdot\left(  3m-1\right)  \cdot\left(  3m\right)  \right)
}{3^{3}m^{3}\left(  3^{m-1}\left(  m-1\right)  !\right)  ^{3}}\nonumber\\
&  =\dfrac{\left(  3\left(  m-1\right)  \right)  !}{\left(  3^{m-1}\left(
m-1\right)  !\right)  ^{3}}\cdot\underbrace{\dfrac{\left(  3m-2\right)
\cdot\left(  3m-1\right)  \cdot\left(  3m\right)  }{3^{3}m^{3}}}%
_{=\dfrac{3m-2}{3m}\cdot\dfrac{3m-1}{3m}}\nonumber\\
&  =\dfrac{\left(  3\left(  m-1\right)  \right)  !}{\left(  3^{m-1}\left(
m-1\right)  !\right)  ^{3}}\cdot\dfrac{3m-2}{3m}\cdot\dfrac{3m-1}{3m}.
\label{sol.bin.-1/2.c.2nd.4}%
\end{align}

On the other hand, Lemma \ref{lem.sol.bin.-1/2.m-1} (applied to $q=-1/3$)
yields%
\begin{equation}
\dbinom{-1/3}{m}=\underbrace{\dfrac{-1/3-m+1}{m}}_{=-\dfrac{3m-2}{3m}}%
\dbinom{-1/3}{m-1}=-\dfrac{3m-2}{3m}\dbinom{-1/3}{m-1}.
\label{sol.bin.-1/2.c.2nd.5a}%
\end{equation}
Also, Lemma \ref{lem.sol.bin.-1/2.m-1} (applied to $q=-2/3$) yields%
\begin{equation}
\dbinom{-2/3}{m}=\underbrace{\dfrac{-2/3-m+1}{m}}_{=-\dfrac{3m-1}{3m}}%
\dbinom{-2/3}{m-1}=-\dfrac{3m-1}{3m}\dbinom{-2/3}{m-1}.
\label{sol.bin.-1/2.c.2nd.5b}%
\end{equation}
Multiplying the two equalities (\ref{sol.bin.-1/2.c.2nd.5a}) and
(\ref{sol.bin.-1/2.c.2nd.5b}), we obtain%
\begin{align*}
\dbinom{-1/3}{m}\dbinom{-2/3}{m}  &  =\left(  -\dfrac{3m-2}{3m}\dbinom
{-1/3}{m-1}\right)  \left(  -\dfrac{3m-1}{3m}\dbinom{-2/3}{m-1}\right) \\
&  =\underbrace{\dbinom{-1/3}{m-1}\dbinom{-2/3}{m-1}}_{\substack{=\dfrac
{\left(  3\left(  m-1\right)  \right)  !}{\left(  3^{m-1}\left(  m-1\right)
!\right)  ^{3}}\\\text{(by (\ref{sol.bin.-1/2.c.2nd.IH}))}}}\cdot\dfrac
{3m-2}{3m}\cdot\dfrac{3m-1}{3m}\\
&  =\dfrac{\left(  3\left(  m-1\right)  \right)  !}{\left(  3^{m-1}\left(
m-1\right)  !\right)  ^{3}}\cdot\dfrac{3m-2}{3m}\cdot\dfrac{3m-1}{3m}%
=\dfrac{\left(  3m\right)  !}{\left(  3^{m}m!\right)  ^{3}}%
\end{align*}
(by (\ref{sol.bin.-1/2.c.2nd.4})). In other words, Exercise \ref{exe.bin.-1/2}
\textbf{(c)} holds for $n=m$. This completes the induction step. Thus, the
induction proof of Exercise \ref{exe.bin.-1/2} \textbf{(c)} is complete.
\end{proof}

\subsubsection{A more general formula}

Parts \textbf{(b)} and \textbf{(c)} of Exercise \ref{exe.bin.-1/2} can be
generalized as follows:

\begin{theorem}
\label{thm.sol.bin.-1/2.gen}Let $h$ be a positive integer. Let $n\in
\mathbb{N}$. Then,%
\[
\prod_{g=1}^{h-1}\dbinom{-g/h}{n}=\left(  \dfrac{-1}{h}\right)  ^{n\left(
h-1\right)  }\cdot\dfrac{\left(  hn\right)  !}{h^{n}n!^{h}}.
\]

\end{theorem}

Exercise \ref{exe.bin.-1/2} \textbf{(b)} follows from Theorem
\ref{thm.sol.bin.-1/2.gen} (applied to $h=2$), after some simple
transformations (using Lemma \ref{lem.sol.bin.-1/2.2nn}). Exercise
\ref{exe.bin.-1/2} \textbf{(c)} follows from Theorem
\ref{thm.sol.bin.-1/2.gen} (applied to $h=3$).

We are going to prove Theorem \ref{thm.sol.bin.-1/2.gen} in a way that is
similar to our second solutions of parts \textbf{(b)} and \textbf{(c)} of
Exercise \ref{exe.bin.-1/2}:

\begin{proof}
[Proof of Theorem \ref{thm.sol.bin.-1/2.gen}.]We shall prove Theorem
\ref{thm.sol.bin.-1/2.gen} by induction on $n$:

\textit{Induction base:} Using $0\left(  h-1\right)  =0$, $h\cdot0=0$ and
$0!=1$, we obtain
\begin{align*}
\left(  \dfrac{-1}{h}\right)  ^{0\left(  h-1\right)  }\cdot\dfrac{\left(
h\cdot0\right)  !}{h^{0}0!^{h}}  &  =\underbrace{\left(  \dfrac{-1}{h}\right)
^{0}}_{=1}\cdot\dfrac{0!}{h^{0}\cdot1^{h}}=\dfrac{0!}{h^{0}\cdot1^{h}}%
=\dfrac{1}{1\cdot1}\\
&  \ \ \ \ \ \ \ \ \ \ \left(  \text{since }0!=1\text{ and }h^{0}=1\text{ and
}1^{h}=1\right) \\
&  =1.
\end{align*}
Comparing this with%
\[
\prod_{g=1}^{h-1}\underbrace{\dbinom{-g/h}{0}}_{\substack{=1\\\text{(by
Proposition \ref{prop.binom.00} \textbf{(a)} (applied to }m=-g/h\text{))}%
}}=\prod_{g=1}^{h-1}1=1,
\]
we obtain $\prod_{g=1}^{h-1}\dbinom{-g/h}{0}=\left(  \dfrac{-1}{h}\right)
^{0\left(  h-1\right)  }\cdot\dfrac{\left(  h\cdot0\right)  !}{h^{0}0!^{h}}$.
In other words, Theorem \ref{thm.sol.bin.-1/2.gen} holds for $n=0$. This
completes the induction base.

\textit{Induction step:} Let $m$ be a positive integer. Assume that Theorem
\ref{thm.sol.bin.-1/2.gen} holds for $n=m-1$. We must now prove that Theorem
\ref{thm.sol.bin.-1/2.gen} holds for $n=m$.

We have $h>0$ (since $h$ is a positive integer) and $m>0$ (since $m$ is a
positive integer). Thus, $hm>0$, so that $hm\neq0$.

We have assumed that Theorem \ref{thm.sol.bin.-1/2.gen} holds for $n=m-1$. In
other words, we have%
\begin{equation}
\prod_{g=1}^{h-1}\dbinom{-g/h}{m-1}=\left(  \dfrac{-1}{h}\right)  ^{\left(
m-1\right)  \left(  h-1\right)  }\cdot\dfrac{\left(  h\left(  m-1\right)
\right)  !}{h^{m-1}\left(  m-1\right)  !^{h}}.
\label{pf.thm.sol.bin.-1/2.gen.IH}%
\end{equation}

We have $m-1\in\mathbb{N}$ (since $m$ is a positive integer) and
$h\in\mathbb{N}$ (since $h$ is a positive integer). Thus, $h\left(
m-1\right)  \in\mathbb{N}$. The definition of $\left(  h\left(  m-1\right)
\right)  !$ thus yields%
\begin{equation}
\left(  h\left(  m-1\right)  \right)  !=1\cdot2\cdot\cdots\cdot\left(
h\left(  m-1\right)  \right)  =\prod_{i=1}^{h\left(  m-1\right)  }i.
\label{pf.thm.sol.bin.-1/2.gen.1}%
\end{equation}

We have $h\left(  m-1\right)  \in\mathbb{N}$ and thus $0\leq h\left(
m-1\right)  $. Also, $hm-h\left(  m-1\right)  =h\geq0$ (since $h\in\mathbb{N}%
$), so that $h\left(  m-1\right)  \leq hm$. The definition of $\left(
hm\right)  !$ yields
\begin{align}
\left(  hm\right)  !  &  =1\cdot2\cdot\cdots\cdot\left(  hm\right)
=\prod_{i=1}^{hm}i=\underbrace{\left(  \prod_{i=1}^{h\left(  m-1\right)
}i\right)  }_{\substack{=\left(  h\left(  m-1\right)  \right)  !\\\text{(by
(\ref{pf.thm.sol.bin.-1/2.gen.1}))}}}\left(  \prod_{i=h\left(  m-1\right)
+1}^{hm}i\right) \nonumber\\
&  \ \ \ \ \ \ \ \ \ \ \left(  \text{since }0\leq h\left(  m-1\right)  \leq
hm\right) \nonumber\\
&  =\left(  h\left(  m-1\right)  \right)  !\cdot\left(  \prod_{i=h\left(
m-1\right)  +1}^{hm}i\right)  . \label{pf.thm.sol.bin.-1/2.gen.1b}%
\end{align}

We have $h\left(  m-1\right)  +1=\left(  hm-h\right)  +1$ and $hm=\left(
hm-h\right)  +h$. Thus,%
\[
\prod_{i=h\left(  m-1\right)  +1}^{hm}i=\prod_{i=\left(  hm-h\right)
+1}^{\left(  hm-h\right)  +h}i=\prod_{g=1}^{h}\left(  g+\left(  hm-h\right)
\right)
\]
(here, we have substituted $g+\left(  hm-h\right)  $ for $i$ in the product).
Thus,%
\begin{align*}
\prod_{i=h\left(  m-1\right)  +1}^{hm}i  &  =\prod_{g=1}^{h}\left(  g+\left(
hm-h\right)  \right)  =\left(  \prod_{g=1}^{h-1}\left(  g+\left(  hm-h\right)
\right)  \right)  \cdot\underbrace{\left(  h+\left(  hm-h\right)  \right)
}_{=hm}\\
&  \ \ \ \ \ \ \ \ \ \ \left(
\begin{array}
[c]{c}%
\text{here, we have split off the factor for }g=h\text{ from}\\
\text{the product (since }h>0\text{)}%
\end{array}
\right) \\
&  =\left(  \prod_{g=1}^{h-1}\left(  g+\left(  hm-h\right)  \right)  \right)
\cdot hm.
\end{align*}
We can divide this equality by $hm$ (because $hm\neq0$), and thus obtain%
\begin{equation}
\dfrac{1}{hm}\prod_{i=h\left(  m-1\right)  +1}^{hm}i=\prod_{g=1}^{h-1}\left(
g+\left(  hm-h\right)  \right)  . \label{pf.thm.sol.bin.-1/2.gen.2a}%
\end{equation}
Now,%
\begin{align}
&  \prod_{g=1}^{h-1}\underbrace{\dfrac{-g/h-m+1}{m}}_{=\left(  -\dfrac{1}%
{h}\right)  \dfrac{1}{m}\left(  g+\left(  hm-h\right)  \right)  }\nonumber\\
&  =\prod_{g=1}^{h-1}\left(  \left(  -\dfrac{1}{h}\right)  \dfrac{1}{m}\left(
g+\left(  hm-h\right)  \right)  \right) \nonumber\\
&  =\underbrace{\left(  \prod_{g=1}^{h-1}\left(  -\dfrac{1}{h}\right)
\right)  }_{=\left(  -\dfrac{1}{h}\right)  ^{h-1}}\underbrace{\left(
\prod_{g=1}^{h-1}\dfrac{1}{m}\right)  }_{\substack{=\left(  \dfrac{1}%
{m}\right)  ^{h-1}\\=\dfrac{1}{m^{h-1}}}}\underbrace{\left(  \prod_{g=1}%
^{h-1}\left(  g+\left(  hm-h\right)  \right)  \right)  }_{\substack{=\dfrac
{1}{hm}\prod_{i=h\left(  m-1\right)  +1}^{hm}i\\\text{(by
(\ref{pf.thm.sol.bin.-1/2.gen.2a}))}}}\nonumber\\
&  =\left(  -\dfrac{1}{h}\right)  ^{h-1}\dfrac{1}{m^{h-1}}\cdot\dfrac{1}%
{hm}\prod_{i=h\left(  m-1\right)  +1}^{hm}i=\left(  -\dfrac{1}{h}\right)
^{h-1}\underbrace{\dfrac{1}{m^{h-1}m}}_{\substack{=\dfrac{1}{m^{h}%
}\\\text{(since }m^{h-1}m=m^{h}\text{)}}}\cdot\dfrac{1}{h}\prod_{i=h\left(
m-1\right)  +1}^{hm}i\nonumber\\
&  =\left(  -\dfrac{1}{h}\right)  ^{h-1}\dfrac{1}{m^{h}}\cdot\dfrac{1}{h}%
\prod_{i=h\left(  m-1\right)  +1}^{hm}i. \label{pf.thm.sol.bin.-1/2.gen.2b}%
\end{align}

We have%
\begin{align}
&  \prod_{g=1}^{h-1}\underbrace{\dbinom{-g/h}{m}}_{\substack{=\dfrac
{-g/h-m+1}{m}\dbinom{-g/h}{m-1}\\\text{(by Lemma \ref{lem.sol.bin.-1/2.m-1}
(applied to }q=-g/h\text{))}}}\nonumber\\
&  =\prod_{g=1}^{h-1}\left(  \dfrac{-g/h-m+1}{m}\dbinom{-g/h}{m-1}\right)
\nonumber\\
&  =\underbrace{\left(  \prod_{g=1}^{h-1}\dfrac{-g/h-m+1}{m}\right)
}_{\substack{=\left(  -\dfrac{1}{h}\right)  ^{h-1}\dfrac{1}{m^{h}}\cdot
\dfrac{1}{h}\prod_{i=h\left(  m-1\right)  +1}^{hm}i\\\text{(by
(\ref{pf.thm.sol.bin.-1/2.gen.2b}))}}}\underbrace{\left(  \prod_{g=1}%
^{h-1}\dbinom{-g/h}{m-1}\right)  }_{\substack{=\left(  \dfrac{-1}{h}\right)
^{\left(  m-1\right)  \left(  h-1\right)  }\cdot\dfrac{\left(  h\left(
m-1\right)  \right)  !}{h^{m-1}\left(  m-1\right)  !^{h}}\\\text{(by
(\ref{pf.thm.sol.bin.-1/2.gen.IH}))}}}\nonumber\\
&  =\left(  -\dfrac{1}{h}\right)  ^{h-1}\dfrac{1}{m^{h}}\cdot\dfrac{1}%
{h}\left(  \prod_{i=h\left(  m-1\right)  +1}^{hm}i\right)  \cdot\left(
\dfrac{-1}{h}\right)  ^{\left(  m-1\right)  \left(  h-1\right)  }\cdot
\dfrac{\left(  h\left(  m-1\right)  \right)  !}{h^{m-1}\left(  m-1\right)
!^{h}}\nonumber\\
&  =\underbrace{\left(  -\dfrac{1}{h}\right)  ^{h-1}\cdot\left(  \dfrac{-1}%
{h}\right)  ^{\left(  m-1\right)  \left(  h-1\right)  }}_{\substack{=\left(
-\dfrac{1}{h}\right)  ^{\left(  h-1\right)  +\left(  m-1\right)  \left(
h-1\right)  }\\=\left(  -\dfrac{1}{h}\right)  ^{m\left(  h-1\right)
}\\\text{(since }\left(  h-1\right)  +\left(  m-1\right)  \left(  h-1\right)
=m\left(  h-1\right)  \text{)}}}\cdot\left(  \prod_{i=h\left(  m-1\right)
+1}^{hm}i\right)  \cdot\dfrac{1}{m^{h}}\cdot\dfrac{1}{h}\dfrac{\left(
h\left(  m-1\right)  \right)  !}{h^{m-1}\left(  m-1\right)  !^{h}}\nonumber\\
&  =\left(  -\dfrac{1}{h}\right)  ^{m\left(  h-1\right)  }\cdot\left(
\prod_{i=h\left(  m-1\right)  +1}^{hm}i\right)  \cdot\dfrac{1}{m^{h}}%
\cdot\dfrac{1}{h}\dfrac{\left(  h\left(  m-1\right)  \right)  !}%
{h^{m-1}\left(  m-1\right)  !^{h}}\nonumber\\
&  =\left(  -\dfrac{1}{h}\right)  ^{m\left(  h-1\right)  }\cdot
\underbrace{\left(  h\left(  m-1\right)  \right)  !\cdot\left(  \prod
_{i=h\left(  m-1\right)  +1}^{hm}i\right)  }_{\substack{=\left(  hm\right)
!\\\text{(by (\ref{pf.thm.sol.bin.-1/2.gen.1b}))}}}\cdot\dfrac{1}{m^{h}}%
\cdot\underbrace{\dfrac{1}{hh^{m-1}\left(  m-1\right)  !^{h}}}%
_{\substack{=\dfrac{1}{h^{m}\left(  m-1\right)  !^{h}}\\\text{(since }%
hh^{m-1}=h^{m}\text{)}}}\nonumber\\
&  =\left(  -\dfrac{1}{h}\right)  ^{m\left(  h-1\right)  }\cdot\left(
hm\right)  !\cdot\dfrac{1}{m^{h}}\cdot\dfrac{1}{h^{m}\left(  m-1\right)
!^{h}}\nonumber\\
&  =\left(  -\dfrac{1}{h}\right)  ^{m\left(  h-1\right)  }\cdot\dfrac{\left(
hm\right)  !}{h^{m}\left(  m-1\right)  !^{h}m^{h}}.
\label{pf.thm.sol.bin.-1/2.gen.5}%
\end{align}

Applying (\ref{eq.n!.rec}) to $n=m$, we obtain $m!=m\cdot\left(  m-1\right)
!$. Thus,%
\[
m!^{h}=\left(  m\cdot\left(  m-1\right)  !\right)  ^{h}=m^{h}\cdot\left(
m-1\right)  !^{h}.
\]
Hence,%
\begin{align*}
\left(  -\dfrac{1}{h}\right)  ^{m\left(  h-1\right)  }\cdot\dfrac{\left(
hm\right)  !}{h^{m}m!^{h}}  &  =\left(  -\dfrac{1}{h}\right)  ^{m\left(
h-1\right)  }\cdot\dfrac{\left(  hm\right)  !}{h^{m}m^{h}\cdot\left(
m-1\right)  !^{h}}\\
&  =\left(  -\dfrac{1}{h}\right)  ^{m\left(  h-1\right)  }\cdot\dfrac{\left(
hm\right)  !}{h^{m}\left(  m-1\right)  !^{h}m^{h}}.
\end{align*}
Comparing this with (\ref{pf.thm.sol.bin.-1/2.gen.5}), we obtain%
\[
\prod_{g=1}^{h-1}\dbinom{-g/h}{m}=\left(  -\dfrac{1}{h}\right)  ^{m\left(
h-1\right)  }\cdot\dfrac{\left(  hm\right)  !}{h^{m}m!^{h}}.
\]
In other words, Theorem \ref{thm.sol.bin.-1/2.gen} holds for $n=m$. This
completes the induction step. Thus, the induction proof of Theorem
\ref{thm.sol.bin.-1/2.gen} is complete.
\end{proof}

\subsection{Solution to Exercise \ref{exe.binom.hockey1}}

Let us start by proving part \textbf{(a)} of Exercise \ref{exe.binom.hockey1}:

\begin{lemma}
\label{lem.sol.multichoose.hock}Let $n\in\mathbb{N}$ and $q\in\mathbb{Q}$.
Then,%
\[
\sum_{r=0}^{n}\dbinom{r+q}{r}=\dbinom{n+q+1}{n}.
\]

\end{lemma}

\begin{proof}
[Proof of Lemma \ref{lem.sol.multichoose.hock}.]Let us forget that $n$ is
fixed. Now, let us prove Lemma \ref{lem.sol.multichoose.hock} by induction on
$n$:

\begin{vershort}
\textit{Induction base:} Proposition \ref{prop.binom.00} \textbf{(a)} (applied
to $m=0+q$) yields $\dbinom{0+q}{0}=1$. Similarly, $\dbinom{0+q+1}{0}=1$.
Comparing this with $\sum_{r=0}^{0}\dbinom{r+q}{r}=\dbinom{0+q}{0}=1$, we
obtain $\sum_{r=0}^{0}\dbinom{r+q}{r}=\dbinom{0+q+1}{0}$. In other words,
Lemma \ref{lem.sol.multichoose.hock} holds for $n=0$. This completes the
induction base.
\end{vershort}

\begin{verlong}
\textit{Induction base:} Proposition \ref{prop.binom.00} \textbf{(a)} (applied
to $0+q+1$ instead of $m$) yields $\dbinom{0+q+1}{0}=1$. On the other hand,
Proposition \ref{prop.binom.00} \textbf{(a)} (applied to $0+q$ instead of $m$)
yields $\dbinom{0+q}{0}=1$. Now, $\sum_{r=0}^{0}\dbinom{r+q}{r}=\dbinom
{0+q}{0}=1=\dbinom{0+q+1}{0}$ (since $\dbinom{0+q+1}{0}=1$). In other words,
Lemma \ref{lem.sol.multichoose.hock} holds for $n=0$. This completes the
induction base.
\end{verlong}

\textit{Induction step:} Let $m$ be a positive integer. Assume that Lemma
\ref{lem.sol.multichoose.hock} holds for $n=m-1$. We must prove that Lemma
\ref{lem.sol.multichoose.hock} holds for $n=m$.

We have assumed that Lemma \ref{lem.sol.multichoose.hock} holds for $n=m-1$.
In other words, we have%
\[
\sum_{r=0}^{m-1}\dbinom{r+q}{r}=\dbinom{\left(  m-1\right)  +q+1}{m-1}.
\]

We have $m\in\left\{  1,2,3,\ldots\right\}  $ (since $m$ is a positive
integer). Thus, Proposition \ref{prop.binom.rec} (applied to $m+q+1$ and $m$
instead of $m$ and $n$) yields
\begin{align}
\dbinom{m+q+1}{m}  &  =\dbinom{\left(  m+q+1\right)  -1}{m}+\dbinom{\left(
m+q+1\right)  -1}{m-1}\nonumber\\
&  =\dbinom{m+q}{m}+\dbinom{m+q}{m-1} \label{pf.lem.sol.multichoose.hock.5}%
\end{align}
(since $\left(  m+q+1\right)  -1=m+q$).

On the other hand, $m\in\left\{  0,1,\ldots,m\right\}  $ (since $m$ is a
positive integer). Thus, we can split off the addend for $r=m$ from the sum
$\sum_{r=0}^{m}\dbinom{r+q}{r}$. We thus obtain%
\begin{align*}
\sum_{r=0}^{m}\dbinom{r+q}{r}  &  =\underbrace{\sum_{r=0}^{m-1}\dbinom{r+q}%
{r}}_{\substack{=\dbinom{\left(  m-1\right)  +q+1}{m-1}\\=\dbinom{m+q}%
{m-1}\\\text{(since }\left(  m-1\right)  +q+1=m+q\text{)}}}+\dbinom{m+q}{m}\\
&  =\dbinom{m+q}{m-1}+\dbinom{m+q}{m}=\dbinom{m+q}{m}+\dbinom{m+q}%
{m-1}=\dbinom{m+q+1}{m}%
\end{align*}
(by (\ref{pf.lem.sol.multichoose.hock.5})). In other words, Lemma
\ref{lem.sol.multichoose.hock} holds for $n=m$. This completes the induction
step. Thus, Lemma \ref{lem.sol.multichoose.hock} is proven by induction.
\end{proof}

Next, let us show the two equalities in Exercise \ref{exe.binom.hockey1}
\textbf{(b)} separately:

\begin{lemma}
\label{lem.sol.binom.hockey1.b1}Let $n\in\left\{  -1,0,1,\ldots\right\}  $ and
$k\in\mathbb{N}$. Then:

\textbf{(a)} We have%
\[
\sum_{i=k}^{n}\dbinom{i}{k}=\dbinom{n+1}{k+1}.
\]

\textbf{(b)} We have%
\[
\sum_{i=0}^{n}\dbinom{i}{k}=\sum_{i=k}^{n}\dbinom{i}{k}.
\]

\end{lemma}

\begin{proof}
[Proof of Lemma \ref{lem.sol.binom.hockey1.b1}.]We have $n\in\left\{
-1,0,1,\ldots\right\}  $, so that $n\geq-1$ and therefore $n+1\geq0$. Hence,
$n+1\in\mathbb{N}$.

\textbf{(a)} If $n<k$, then Lemma \ref{lem.sol.binom.hockey1.b1} \textbf{(a)}
holds\footnote{\textit{Proof.} Assume that $n<k$. We must show that Lemma
\ref{lem.sol.binom.hockey1.b1} \textbf{(a)} holds.
\par
Recall that $n+1\in\mathbb{N}$. Also, $k+1>k\geq0$ (since $k\in\mathbb{N}$)
and thus $k+1\in\mathbb{N}$. Finally, $\underbrace{n}_{<k}+1<k+1$. Hence,
Proposition \ref{prop.binom.0} (applied to $n+1$ and $k+1$ instead of $m$ and
$n$) yields $\dbinom{n+1}{k+1}=0$. Comparing this with%
\begin{align*}
\sum_{i=k}^{n}\dbinom{i}{k}  &  =\left(  \text{empty sum}\right)
\ \ \ \ \ \ \ \ \ \ \left(  \text{since }n<k\right) \\
&  =0,
\end{align*}
we obtain $\sum_{i=k}^{n}\dbinom{i}{k}=\dbinom{n+1}{k+1}$. In other words,
Lemma \ref{lem.sol.binom.hockey1.b1} \textbf{(a)} holds. Thus, we have shown
that if $n<k$, then Lemma \ref{lem.sol.binom.hockey1.b1} \textbf{(a)} holds.}.
Hence, for the rest of this proof of Lemma \ref{lem.sol.binom.hockey1.b1}
\textbf{(a)}, we can WLOG assume that we don't have $n<k$. Assume this.

We have $n\geq k$ (since we don't have $n<k$). Hence, $n\geq k\geq0$ (since
$k\in\mathbb{N}$), so that $n\in\mathbb{N}$. Also, from $n\geq k$, we obtain
$n-k\geq0$, so that $n-k\in\mathbb{N}$.

Lemma \ref{lem.sol.multichoose.hock} (applied to $n-k$ and $k$ instead of $n$
and $q$) yields%
\begin{equation}
\sum_{r=0}^{n-k}\dbinom{r+k}{r}=\dbinom{\left(  n-k\right)  +k+1}{n-k}%
=\dbinom{n+1}{n-k} \label{pf.lem.sol.binom.hockey1.b1.1}%
\end{equation}
(since $\left(  n-k\right)  +k+1=n+1$).

We have $n+1\in\mathbb{N}$ and $k+1\in\mathbb{N}$ (since $k\in\mathbb{N}$).
Also, $\underbrace{n}_{\geq k}+1\geq k+1$. Thus, Proposition
\ref{prop.binom.symm} (applied to $n+1$ and $k+1$ instead of $n$ and $k$)
yields
\[
\dbinom{n+1}{k+1}=\dbinom{n+1}{\left(  n+1\right)  -\left(  k+1\right)
}=\dbinom{n+1}{n-k}%
\]
(since $\left(  n+1\right)  -\left(  k+1\right)  =n-k$). Comparing this
equality with (\ref{pf.lem.sol.binom.hockey1.b1.1}), we obtain%
\begin{equation}
\dbinom{n+1}{k+1}=\sum_{r=0}^{n-k}\dbinom{r+k}{r}=\sum_{i=k}^{n}%
\dbinom{\left(  i-k\right)  +k}{i-k} \label{pf.lem.sol.binom.hockey1.b1.4}%
\end{equation}
(here, we have substituted $i-k$ for $r$ in the sum).

But each $i\in\left\{  k,k+1,\ldots,n\right\}  $ satisfies
\begin{equation}
\dbinom{\left(  i-k\right)  +k}{i-k}=\dbinom{i}{k}.
\label{pf.lem.sol.binom.hockey1.b1.2}%
\end{equation}

[\textit{Proof of (\ref{pf.lem.sol.binom.hockey1.b1.2}):} Let $i\in\left\{
k,k+1,\ldots,n\right\}  $. Thus, $i\geq k$ and $i\in\left\{  k,k+1,\ldots
,n\right\}  \subseteq\mathbb{N}$. Hence, Proposition \ref{prop.binom.symm}
(applied to $i$ and $k$ instead of $m$ and $n$) yields $\dbinom{i}{k}%
=\dbinom{i}{i-k}$. Now, $\left(  i-k\right)  +k=i$, so that $\dbinom{\left(
i-k\right)  +k}{i-k}=\dbinom{i}{i-k}=\dbinom{i}{k}$. This proves
(\ref{pf.lem.sol.binom.hockey1.b1.2}).]

Now, (\ref{pf.lem.sol.binom.hockey1.b1.4}) becomes%
\[
\dbinom{n+1}{k+1}=\sum_{i=k}^{n}\underbrace{\dbinom{\left(  i-k\right)
+k}{i-k}}_{\substack{=\dbinom{i}{k}\\\text{(by
(\ref{pf.lem.sol.binom.hockey1.b1.2}))}}}=\sum_{i=k}^{n}\dbinom{i}{k}.
\]
In other words, $\sum_{i=k}^{n}\dbinom{i}{k}=\dbinom{n+1}{k+1}$. This proves
Lemma \ref{lem.sol.binom.hockey1.b1} \textbf{(a)}.

\textbf{(b)} If $n<k$, then Lemma \ref{lem.sol.binom.hockey1.b1} \textbf{(b)}
holds\footnote{\textit{Proof.} Assume that $n<k$. We must show that Lemma
\ref{lem.sol.binom.hockey1.b1} \textbf{(b)} holds.
\par
Let $i\in\left\{  0,1,\ldots,n\right\}  $. Thus, $i\leq n<k$ and $i\in\left\{
0,1,\ldots,n\right\}  \subseteq\mathbb{N}$. Hence, Proposition
\ref{prop.binom.0} (applied to $i$ and $k$ instead of $m$ and $n$) yields
$\dbinom{i}{k}=0$.
\par
Now, forget that we fixed $i$. We thus have proven the equality $\dbinom{i}%
{k}=0$ for each $i\in\left\{  0,1,\ldots,n\right\}  $. Adding up these
equalities for all $i\in\left\{  0,1,\ldots,n\right\}  $, we obtain
$\sum_{i=0}^{n}\dbinom{i}{k}=\sum_{i=0}^{n}0=0$. Comparing this with%
\begin{align*}
\sum_{i=k}^{n}\dbinom{i}{k}  &  =\left(  \text{empty sum}\right)
\ \ \ \ \ \ \ \ \ \ \left(  \text{since }n<k\right) \\
&  =0,
\end{align*}
we obtain $\sum_{i=0}^{n}\dbinom{i}{k}=\sum_{i=k}^{n}\dbinom{i}{k}$. In other
words, Lemma \ref{lem.sol.binom.hockey1.b1} \textbf{(b)} holds. Thus, we have
shown that if $n<k$, then Lemma \ref{lem.sol.binom.hockey1.b1} \textbf{(b)}
holds.}. Hence, for the rest of this proof, we can WLOG assume that we don't
have $n<k$. Assume this.

Let $i\in\left\{  0,1,\ldots,k-1\right\}  $. Thus, $i\leq k-1<k$ and
$i\in\left\{  0,1,\ldots,k-1\right\}  \subseteq\mathbb{N}$. Hence, Proposition
\ref{prop.binom.0} (applied to $i$ and $k$ instead of $m$ and $n$) yields
$\dbinom{i}{k}=0$.

Now, forget that we fixed $i$. We thus have proven the equality $\dbinom{i}%
{k}=0$ for each $i\in\left\{  0,1,\ldots,k-1\right\}  $. Adding up these
equalities for all $i\in\left\{  0,1,\ldots,k-1\right\}  $, we obtain
$\sum_{i=0}^{k-1}\dbinom{i}{k}=\sum_{i=0}^{k-1}0=0$.

From $k\in\mathbb{N}$, we obtain $k\geq0$, so that $0\leq k\leq n$ (since
$n\geq k$ (because we don't have $n<k$)). Hence, we can split the sum
$\sum_{i=0}^{n}\dbinom{i}{k}$ at $i=k$. We thus obtain%
\[
\sum_{i=0}^{n}\dbinom{i}{k}=\underbrace{\sum_{i=0}^{k-1}\dbinom{i}{k}}%
_{=0}+\sum_{i=k}^{n}\dbinom{i}{k}=\sum_{i=k}^{n}\dbinom{i}{k}.
\]
This proves Lemma \ref{lem.sol.binom.hockey1.b1} \textbf{(b)}.
\end{proof}

\begin{proof}
[Solution to Exercise \ref{exe.binom.hockey1}.]\textbf{(a)} Let $n\in
\mathbb{N}$ and $q\in\mathbb{Q}$. Then, Lemma \ref{lem.sol.multichoose.hock}
yields%
\[
\sum_{r=0}^{n}\dbinom{r+q}{r}=\dbinom{n+q+1}{n}.
\]
This solves Exercise \ref{exe.binom.hockey1} \textbf{(a)}.

\textbf{(b)} Let $n\in\left\{  -1,0,1,\ldots\right\}  $ and $k\in\mathbb{N}$.
Lemma \ref{lem.sol.binom.hockey1.b1} \textbf{(b)} yields%
\[
\sum_{i=0}^{n}\dbinom{i}{k}=\sum_{i=k}^{n}\dbinom{i}{k}=\dbinom{n+1}{k+1}%
\]
(by Lemma \ref{lem.sol.binom.hockey1.b1} \textbf{(a)}). This solves Exercise
\ref{exe.binom.hockey1} \textbf{(b)}.
\end{proof}

\subsection{Solution to Exercise \ref{exe.prop.binom.subsets}}

\begin{vershort}
Exercise \ref{exe.prop.binom.subsets} asks us to prove Proposition
\ref{prop.binom.subsets}.

First, let us handle two trivial cases of Proposition \ref{prop.binom.subsets}:

\begin{lemma}
\label{lem.sol.prop.binom.subsets.n=0}Proposition \ref{prop.binom.subsets}
holds in the case when $n=0$.
\end{lemma}

\begin{proof}
[Proof of Lemma \ref{lem.sol.prop.binom.subsets.n=0}.]Let $m\in\mathbb{N}$,
and let $S$ be an $m$-element set. Then, the set $S$ has exactly one
$0$-element subset (namely, the empty set $\varnothing$). Thus, the number of
all $0$-element subsets of $S$ is $1$.

On the other hand, $\dbinom{m}{0}=1$ (by (\ref{eq.binom.00})). Thus,
$\dbinom{m}{0}$ is the number of all $0$-element subsets of $S$ (since the
number of all $0$-element subsets of $S$ is $1$).

Now, forget that we fixed $m$ and $S$. We thus have shown that if
$m\in\mathbb{N}$ and if $S$ is an $m$-element set, then $\dbinom{m}{0}$ is the
number of all $0$-element subsets of $S$. In other words, Proposition
\ref{prop.binom.subsets} holds in the case when $n=0$. This proves Lemma
\ref{lem.sol.prop.binom.subsets.n=0}.
\end{proof}

\begin{lemma}
\label{lem.sol.prop.binom.subsets.m=0}Proposition \ref{prop.binom.subsets}
holds in the case when $m=0$.
\end{lemma}

\begin{proof}
[Proof of Lemma \ref{lem.sol.prop.binom.subsets.m=0}.]Let $m$, $n$ and $S$ be
as in Proposition \ref{prop.binom.subsets}. Assume that $m=0$. We then must
prove that Proposition \ref{prop.binom.subsets} holds for these $m$, $n$ and
$S$.

If $n=0$, then this follows from Lemma \ref{lem.sol.prop.binom.subsets.n=0}.
Hence, for the rest of this proof, we can WLOG assume that $n\neq0$. Assume this.

From $n\neq0$, we conclude that $n$ is a positive integer. Hence, $n>0=m$, so
that $m<n$. Thus, (\ref{eq.binom.0}) yields $\dbinom{m}{n}=0$.

But $S$ is an $m$-element set, i.e., a $0$-element set (since $m=0$). In other
words, $S=\varnothing$. Hence, $S$ has no $n$-element
subsets\footnote{\textit{Proof.} Assume the contrary. Thus, $S$ has an
$n$-element subset. Let $Q$ be such a subset. Thus, $Q\subseteq S$ and
$\left\vert Q\right\vert =n$. Since $Q\subseteq S=\varnothing$, we have
$Q=\varnothing$ and thus $\left\vert Q\right\vert =0$. This contradicts
$\left\vert Q\right\vert =n\neq0$. This contradiction shows that our
assumption was wrong, qed.}. In other words, the number of all $n$-element
subsets of $S$ is $0$. Since $\dbinom{m}{n}$ is also $0$, this shows that
$\dbinom{m}{n}$ is the number of all $n$-element subsets of $S$. Hence,
Proposition \ref{prop.binom.subsets} holds for our $m$, $n$ and $S$. This
proves Lemma \ref{lem.sol.prop.binom.subsets.m=0}.
\end{proof}

\begin{proof}
[Proof of Proposition \ref{prop.binom.subsets}.]We shall prove Proposition
\ref{prop.binom.subsets} by induction over $m$:

\textit{Induction base:} Lemma \ref{lem.sol.prop.binom.subsets.m=0} shows that
Proposition \ref{prop.binom.subsets} holds in the case when $m=0$. This
completes the induction base.

\textit{Induction step:} Let $M$ be a positive integer. Assume that
Proposition \ref{prop.binom.subsets} holds for $m=M-1$. We now must prove that
Proposition \ref{prop.binom.subsets} holds for $m=M$.

We have assumed that Proposition \ref{prop.binom.subsets} holds for $m=M-1$.
In other words, if $n\in\mathbb{N}$, and if $S$ is an $\left(  M-1\right)
$-element set, then%
\begin{equation}
\dbinom{M-1}{n}\text{ is the number of all }n\text{-element subsets of }S.
\label{pf.prop.binom.subsets.ihyp}%
\end{equation}

Now, let $n\in\mathbb{N}$, and let $S$ be an $M$-element set. We shall show
that
\begin{equation}
\dbinom{M}{n}\text{ is the number of all }n\text{-element subsets of }S.
\label{pf.prop.binom.subsets.igoal}%
\end{equation}

[\textit{Proof of (\ref{pf.prop.binom.subsets.igoal}):} If $n=0$, then
(\ref{pf.prop.binom.subsets.igoal}) follows from Lemma
\ref{lem.sol.prop.binom.subsets.n=0}\footnote{\textit{Proof.} Lemma
\ref{lem.sol.prop.binom.subsets.n=0} yields that Proposition
\ref{prop.binom.subsets} holds for $n=0$. Hence, we can apply Proposition
\ref{prop.binom.subsets} to $n=0$ and $m=M$. We thus obtain that $\dbinom
{M}{0}$ is the number of all $0$-element subsets of $S$. If $n=0$, then this
rewrites as follows: $\dbinom{M}{n}$ is the number of all $n$-element subsets
of $S$. Hence, if $n=0$, then (\ref{pf.prop.binom.subsets.igoal}) holds.}.
Thus, for the rest of this proof of (\ref{pf.prop.binom.subsets.igoal}), we
can WLOG assume that $n\neq0$. Assume this.

Now, $n$ is a positive integer (since $n\in\mathbb{N}$ and $n\neq0$); thus,
$n-1\in\mathbb{N}$.

Now, $S$ is an $M$-element set; thus, $\left\vert S\right\vert =M>0$. Hence,
the set $S$ is nonempty. In other words, there exists an $s\in S$. Fix such an
$s$. Then, $\left\vert S\setminus\left\{  s\right\}  \right\vert
=\underbrace{\left\vert S\right\vert }_{=M}-1=M-1$. In other words,
$S\setminus\left\{  s\right\}  $ is an $\left(  M-1\right)  $-element set.
Hence, (\ref{pf.prop.binom.subsets.ihyp}) (applied to $S\setminus\left\{
s\right\}  $ instead of $S$) yields that%
\begin{equation}
\dbinom{M-1}{n}\text{ is the number of all }n\text{-element subsets of
}S\setminus\left\{  s\right\}  . \label{pf.prop.binom.subsets.ihyp.1}%
\end{equation}
Also, (\ref{pf.prop.binom.subsets.ihyp}) (applied to $n-1$ and $S\setminus
\left\{  s\right\}  $ instead of $n$ and $S$) yields that%
\begin{equation}
\dbinom{M-1}{n-1}\text{ is the number of all }\left(  n-1\right)
\text{-element subsets of }S\setminus\left\{  s\right\}  .
\label{pf.prop.binom.subsets.ihyp.2}%
\end{equation}

Now, the $n$-element subsets of $S$ can be classified into two types: the ones
that contain $s$, and the ones that don't. We can count them separately:

\begin{itemize}
\item The $n$-element subsets of $S$ that contain $s$ are in bijection with
the $\left(  n-1\right)  $-element subsets of $S\setminus\left\{  s\right\}
$. More precisely: To each $\left(  n-1\right)  $-element subset $U$ of
$S\setminus\left\{  s\right\}  $, we can assign a unique $n$-element subset of
$S$ that contains $s$ (namely, $U\cup\left\{  s\right\}  $); and this
assignment is bijective (i.e., each $n$-element subset of $S$ that contains
$s$ gets assigned to exactly one $\left(  n-1\right)  $-element subsets of
$S\setminus\left\{  s\right\}  $).\ \ \ \ \footnote{Let me restate this in
even more formal terms:
\par
Let $\mathbf{A}$ be the set of all $\left(  n-1\right)  $-element subsets of
$S\setminus\left\{  s\right\}  $. Let $\mathbf{B}$ be the set of all
$n$-element subsets of $S$ that contain $s$. Then, the map%
\[
\mathbf{A}\rightarrow\mathbf{B},\ \ \ \ \ \ \ \ \ \ U\mapsto U\cup\left\{
s\right\}
\]
is well-defined and bijective.
\par
(If you want to prove this formally, you need to prove two statements:
\par
\begin{enumerate}
\item The map $\mathbf{A}\rightarrow\mathbf{B},\ U\mapsto U\cup\left\{
s\right\}  $ is well-defined (i.e., we have $U\cup\left\{  s\right\}
\in\mathbf{B}$ for each $U\in\mathbf{A}$).
\par
\item This map is bijective.
\end{enumerate}
\par
Proving the first statement is straightforward. The best way to prove the
second statement is to show that the map $\mathbf{A}\rightarrow\mathbf{B}%
,\ U\mapsto U\cup\left\{  s\right\}  $ has an inverse -- namely, the map
$\mathbf{B}\rightarrow\mathbf{A},\ V\mapsto V\setminus\left\{  s\right\}  $.
Of course, you would also have to show that this latter map is well-defined,
too.)} Hence,
\begin{align}
&  \left(  \text{the number of all }n\text{-element subsets of }S\text{ that
contain }s\right) \nonumber\\
&  =\left(  \text{the number of all }\left(  n-1\right)  \text{-element
subsets of }S\setminus\left\{  s\right\}  \right) \nonumber\\
&  =\dbinom{M-1}{n-1}\ \ \ \ \ \ \ \ \ \ \left(  \text{by
(\ref{pf.prop.binom.subsets.ihyp.2})}\right)  .
\label{pf.prop.binom.subsets.i1}%
\end{align}

\item The $n$-element subsets of $S$ that don't contain $s$ are precisely the
$n$-element subsets of $S\setminus\left\{  s\right\}  $ (because the subsets
of $S$ that don't contain $s$ are precisely the subsets of $S\setminus\left\{
s\right\}  $). Hence,%
\begin{align}
&  \left(  \text{the number of all }n\text{-element subsets of }S\text{ that
don't contain }s\right) \nonumber\\
&  =\left(  \text{the number of all }n\text{-element subsets of }%
S\setminus\left\{  s\right\}  \right) \nonumber\\
&  =\dbinom{M-1}{n}\ \ \ \ \ \ \ \ \ \ \left(  \text{by
(\ref{pf.prop.binom.subsets.ihyp.1})}\right)  .
\label{pf.prop.binom.subsets.i2}%
\end{align}

\end{itemize}

Now, every $n$-element subset of $S$ either contains $s$ or does not. Hence,%
\begin{align*}
&  \left(  \text{the number of all }n\text{-element subsets of }S\right) \\
&  =\underbrace{\left(  \text{the number of all }n\text{-element subsets of
}S\text{ that contain }s\right)  }_{\substack{=\dbinom{M-1}{n-1}\\\text{(by
(\ref{pf.prop.binom.subsets.i1}))}}}\\
&  \ \ \ \ \ \ \ \ \ \ +\underbrace{\left(  \text{the number of all
}n\text{-element subsets of }S\text{ that don't contain }s\right)
}_{\substack{=\dbinom{M-1}{n}\\\text{(by (\ref{pf.prop.binom.subsets.i2}))}%
}}\\
&  =\dbinom{M-1}{n-1}+\dbinom{M-1}{n}.
\end{align*}
Compared with
\[
\dbinom{M}{n}=\dbinom{M-1}{n-1}+\dbinom{M-1}{n}\ \ \ \ \ \ \ \ \ \ \left(
\text{by (\ref{eq.binom.rec.m}), applied to }m=M\right)  ,
\]
this yields $\left(  \text{the number of all }n\text{-element subsets of
}S\right)  =\dbinom{M}{n}$. Hence, (\ref{pf.prop.binom.subsets.igoal}) is proven.]

Now, forget that we fixed $n$ and $S$. We thus have shown that every
$n\in\mathbb{N}$ and every $M$-element set $S$ satisfy
(\ref{pf.prop.binom.subsets.igoal}). In other words, Proposition
\ref{prop.binom.subsets} holds for $m=M$. This completes the induction step.
Thus, Proposition \ref{prop.binom.subsets} is proven by induction.
\end{proof}
\end{vershort}

\begin{verlong}
In order to solve Exercise \ref{exe.prop.binom.subsets}, we need to prove
Proposition \ref{prop.binom.subsets}. We shall do this in detail; but first,
let us introduce a notation:

\begin{definition}
\label{def.sol.prop.binom.subsets.Pm}Let $S$ be a set. Let $m\in\mathbb{N}$.
Then, $\mathcal{P}_{m}\left(  S\right)  $ will denote the set of all
$m$-element subsets of $S$.
\end{definition}

\begin{proposition}
\label{prop.sol.prop.binom.subsets.Pm.lem}\textbf{(a)} We have $\mathcal{P}%
_{0}\left(  S\right)  =\left\{  \varnothing\right\}  $ for every set $S$.

\textbf{(b)} For every positive integer $m$, we have $\mathcal{P}_{m}\left(
\varnothing\right)  =\varnothing$.

\textbf{(c)} Let $S$ be a set. Let $s\in S$. Let $m$ be a positive integer.
Then, $\mathcal{P}_{m}\left(  S\setminus\left\{  s\right\}  \right)
\subseteq\mathcal{P}_{m}\left(  S\right)  $. Furthermore, the map%
\begin{align*}
\mathcal{P}_{m-1}\left(  S\setminus\left\{  s\right\}  \right)   &
\rightarrow\mathcal{P}_{m}\left(  S\right)  \setminus\mathcal{P}_{m}\left(
S\setminus\left\{  s\right\}  \right)  ,\\
U  &  \mapsto U\cup\left\{  s\right\}
\end{align*}
is well-defined and a bijection.
\end{proposition}

Let us briefly explain what Proposition
\ref{prop.sol.prop.binom.subsets.Pm.lem} says. Proposition
\ref{prop.sol.prop.binom.subsets.Pm.lem} \textbf{(a)} says that the only
$0$-element subset of any set $S$ is $\varnothing$. Proposition
\ref{prop.sol.prop.binom.subsets.Pm.lem} \textbf{(b)} says that the empty set
$\varnothing$ has no $m$-element subsets when $m$ is a positive integer.
Proposition \ref{prop.sol.prop.binom.subsets.Pm.lem} \textbf{(c)} says that if
$s$ is an element of a set $S$, and if $m$ is a positive integer, then:

\begin{itemize}
\item all $m$-element subsets of $S\setminus\left\{  s\right\}  $ are
$m$-element subsets of $S$ as well;

\item the $m$-element subsets of $S$ which are \textbf{not} $m$-element
subsets of $S\setminus\left\{  s\right\}  $ are in bijection with the $\left(
m-1\right)  $-element subsets of $S\setminus\left\{  s\right\}  $; this
bijection sends an $\left(  m-1\right)  $-element subset $U$ of $S\setminus
\left\{  s\right\}  $ to the $m$-element subset $U\cup\left\{  s\right\}  $ of
$S$.
\end{itemize}

This restatement should make Proposition
\ref{prop.sol.prop.binom.subsets.Pm.lem} obvious (or, at least, intuitively
clear). For the sake of completeness, let me also give a formal proof:

\begin{proof}
[Proof of Proposition \ref{prop.sol.prop.binom.subsets.Pm.lem}.]\textbf{(a)}
Let $S$ be a set. The set $\varnothing$ is clearly a $0$-element subset of
$S$, and furthermore is the only $0$-element subset of $S$ (since
$\varnothing$ is the only $0$-element set). Now, the definition of
$\mathcal{P}_{0}\left(  S\right)  $ shows that $\mathcal{P}_{0}\left(
S\right)  $ is the set of all $0$-element subsets of $S$. Thus,%
\[
\mathcal{P}_{0}\left(  S\right)  =\left(  \text{the set of all }%
0\text{-element subsets of }S\right)  =\left\{  \varnothing\right\}
\]
(since $\varnothing$ is the only $0$-element subset of $S$). This proves
Proposition \ref{prop.sol.prop.binom.subsets.Pm.lem} \textbf{(a)}.

\textbf{(b)} Let $m$ be a positive integer. Then, there exist no $m$-element
subsets of $\varnothing$\ \ \ \ \footnote{\textit{Proof.} Let $U$ be an
$m$-element subset of $\varnothing$. Thus, $U$ is an $m$-element set; hence,
$\left\vert U\right\vert =m>0$ (since $m$ is positive). But $U$ is a subset of
$\varnothing$; hence, $\left\vert U\right\vert \leq\left\vert \varnothing
\right\vert =0$. This contradicts $\left\vert U\right\vert >0$.
\par
Now, forget that we fixed $U$. We thus have found a contradiction for every
$m$-element subset $U$ of $\varnothing$. Hence, there exist no $m$-element
subsets of $\varnothing$. Qed.}. Now, the definition of $\mathcal{P}%
_{m}\left(  \varnothing\right)  $ shows that $\mathcal{P}_{m}\left(
\varnothing\right)  $ is the set of all $m$-element subsets of $\varnothing$.
Hence,
\[
\mathcal{P}_{m}\left(  \varnothing\right)  =\left(  \text{the set of all
}m\text{-element subsets of }\varnothing\right)  =\varnothing
\]
(since there exist no $m$-element subsets of $\varnothing$). This proves
Proposition \ref{prop.sol.prop.binom.subsets.Pm.lem} \textbf{(b)}.

\textbf{(c)} We have $m-1\in\mathbb{N}$ (since $m$ is a positive integer). We
know that $\mathcal{P}_{m}\left(  S\setminus\left\{  s\right\}  \right)  $ is
the set of all $m$-element subsets of $S\setminus\left\{  s\right\}  $ (by the
definition of $\mathcal{P}_{m}\left(  S\setminus\left\{  s\right\}  \right)
$). Also, $\mathcal{P}_{m}\left(  S\right)  $ is the set of all $m$-element
subsets of $S$ (by the definition of $\mathcal{P}_{m}\left(  S\right)  $).
Finally, $\mathcal{P}_{m-1}\left(  S\setminus\left\{  s\right\}  \right)  $ is
the set of all $\left(  m-1\right)  $-element subsets of $S\setminus\left\{
s\right\}  $ (by the definition of $\mathcal{P}_{m-1}\left(  S\setminus
\left\{  s\right\}  \right)  $).

For every $U\in\mathcal{P}_{m}\left(  S\setminus\left\{  s\right\}  \right)
$, we have $U\in\mathcal{P}_{m}\left(  S\right)  $%
\ \ \ \ \footnote{\textit{Proof.} Let $U\in\mathcal{P}_{m}\left(
S\setminus\left\{  s\right\}  \right)  $. We must show that $U\in
\mathcal{P}_{m}\left(  S\right)  $.
\par
We have $U\in\mathcal{P}_{m}\left(  S\setminus\left\{  s\right\}  \right)  $.
In other words, $U$ is an $m$-element subset of $S\setminus\left\{  s\right\}
$ (since $\mathcal{P}_{m}\left(  S\setminus\left\{  s\right\}  \right)  $ is
the set of all $m$-element subsets of $S\setminus\left\{  s\right\}  $). Thus,
$U$ is an $m$-element set; hence, $\left\vert U\right\vert =m$. Also, $U$ is a
subset of $S\setminus\left\{  s\right\}  $; thus, $U\subseteq S\setminus
\left\{  s\right\}  \subseteq S$. Therefore, $U$ is a subset of $S$. Hence,
$U$ is an $m$-element subset of $S$. In other words, $U\in\mathcal{P}%
_{m}\left(  S\right)  $ (since $\mathcal{P}_{m}\left(  S\right)  $ is the set
of all $m$-element subsets of $S$). Qed.}. In other words, $\mathcal{P}%
_{m}\left(  S\setminus\left\{  s\right\}  \right)  \subseteq\mathcal{P}%
_{m}\left(  S\right)  $.

For every $U\in\mathcal{P}_{m-1}\left(  S\setminus\left\{  s\right\}  \right)
$, we have $U\cup\left\{  s\right\}  \in\mathcal{P}_{m}\left(  S\right)
\setminus\mathcal{P}_{m}\left(  S\setminus\left\{  s\right\}  \right)
$\ \ \ \ \footnote{\textit{Proof.} Let $U\in\mathcal{P}_{m-1}\left(
S\setminus\left\{  s\right\}  \right)  $. We must prove that $U\cup\left\{
s\right\}  \in\mathcal{P}_{m}\left(  S\right)  \setminus\mathcal{P}_{m}\left(
S\setminus\left\{  s\right\}  \right)  $.
\par
Let $V=U\cup\left\{  s\right\}  $.
\par
We have $U\in\mathcal{P}_{m-1}\left(  S\setminus\left\{  s\right\}  \right)
$. In other words, $U$ is an $\left(  m-1\right)  $-element subset of
$S\setminus\left\{  s\right\}  $ (since $\mathcal{P}_{m-1}\left(
S\setminus\left\{  s\right\}  \right)  $ is the set of all $\left(
m-1\right)  $-element subsets of $S\setminus\left\{  s\right\}  $). Thus, $U$
is an $\left(  m-1\right)  $-element set; in other words, $\left\vert
U\right\vert =m-1$. Furthermore, $U$ is a subset of $S\setminus\left\{
s\right\}  $; hence, $U\subseteq S\setminus\left\{  s\right\}  $. Now,
\[
V=\underbrace{U}_{\subseteq S\setminus\left\{  s\right\}  \subseteq S}%
\cup\underbrace{\left\{  s\right\}  }_{\substack{\subseteq S\\\text{(since
}s\in S\text{)}}}\subseteq S\cup S=S.
\]
In other words, $V$ is a subset of $S$.
\par
We have $s\in\left\{  s\right\}  $ and thus $s\notin S\setminus\left\{
s\right\}  $. If we had $s\in U$, then we would have $s\in U\subseteq
S\setminus\left\{  s\right\}  $, which would contradict $s\notin
S\setminus\left\{  s\right\}  $. Thus, we cannot have $s\in U$. In other
words, we have $s\notin U$. Hence, $\left\vert U\cup\left\{  s\right\}
\right\vert =\underbrace{\left\vert U\right\vert }_{=m-1}+1=\left(
m-1\right)  +1=m$. Thus, $\left\vert \underbrace{V}_{=U\cup\left\{  s\right\}
}\right\vert =\left\vert U\cup\left\{  s\right\}  \right\vert =m$. In other
words, $V$ is an $m$-element set. Thus, $V$ is an $m$-element subset of $S$.
In other words, $V\in\mathcal{P}_{m}\left(  S\right)  $ (since $\mathcal{P}%
_{m}\left(  S\right)  $ is the set of all $m$-element subsets of $S$).
\par
Now, we shall prove that $V\notin\mathcal{P}_{m}\left(  S\setminus\left\{
s\right\}  \right)  $. Indeed, assume the contrary (for the sake of
contradiction). Thus, $V\in\mathcal{P}_{m}\left(  S\setminus\left\{
s\right\}  \right)  $. In other words, $V$ is an $m$-element subset of
$S\setminus\left\{  s\right\}  $ (since $\mathcal{P}_{m}\left(  S\setminus
\left\{  s\right\}  \right)  $ is the set of all $m$-element subsets of
$S\setminus\left\{  s\right\}  $). Hence, $V$ is a subset of $S\setminus
\left\{  s\right\}  $. In other words, $V\subseteq S\setminus\left\{
s\right\}  $. Thus, $s\in\left\{  s\right\}  \subseteq U\cup\left\{
s\right\}  =V\subseteq S\setminus\left\{  s\right\}  $, which contradicts
$s\notin S\setminus\left\{  s\right\}  $. This contradiction proves that our
assumption was false. Hence, $V\notin\mathcal{P}_{m}\left(  S\setminus\left\{
s\right\}  \right)  $ is proven.
\par
Combining $V\in\mathcal{P}_{m}\left(  S\right)  $ with $V\notin\mathcal{P}%
_{m}\left(  S\setminus\left\{  s\right\}  \right)  $, we obtain $V\in
\mathcal{P}_{m}\left(  S\right)  \setminus\mathcal{P}_{m}\left(
S\setminus\left\{  s\right\}  \right)  $. Thus, $U\cup\left\{  s\right\}
=V\in\mathcal{P}_{m}\left(  S\right)  \setminus\mathcal{P}_{m}\left(
S\setminus\left\{  s\right\}  \right)  $, qed.}. Hence, we can define a map%
\[
\alpha:\mathcal{P}_{m-1}\left(  S\setminus\left\{  s\right\}  \right)
\rightarrow\mathcal{P}_{m}\left(  S\right)  \setminus\mathcal{P}_{m}\left(
S\setminus\left\{  s\right\}  \right)
\]
by
\begin{equation}
\left(  \alpha\left(  U\right)  =U\cup\left\{  s\right\}
\ \ \ \ \ \ \ \ \ \ \text{for every }U\in\mathcal{P}_{m-1}\left(
S\setminus\left\{  s\right\}  \right)  \right)  .
\label{pf.prop.sol.prop.binom.subsets.Pm.lem.c.alpha()=}%
\end{equation}
Consider this map $\alpha$.

For every $V\in\mathcal{P}_{m}\left(  S\right)  \setminus\mathcal{P}%
_{m}\left(  S\setminus\left\{  s\right\}  \right)  $, we have $V\setminus
\left\{  s\right\}  \in\mathcal{P}_{m-1}\left(  S\setminus\left\{  s\right\}
\right)  $\ \ \ \ \footnote{\textit{Proof.} Let $V\in\mathcal{P}_{m}\left(
S\right)  \setminus\mathcal{P}_{m}\left(  S\setminus\left\{  s\right\}
\right)  $. We must prove that $V\setminus\left\{  s\right\}  \in
\mathcal{P}_{m-1}\left(  S\setminus\left\{  s\right\}  \right)  $.
\par
Let $W=V\setminus\left\{  s\right\}  $.
\par
We have $V\notin\mathcal{P}_{m}\left(  S\setminus\left\{  s\right\}  \right)
$ (since $V\in\mathcal{P}_{m}\left(  S\right)  \setminus\mathcal{P}_{m}\left(
S\setminus\left\{  s\right\}  \right)  $).
\par
We have $V\in\mathcal{P}_{m}\left(  S\right)  \setminus\mathcal{P}_{m}\left(
S\setminus\left\{  s\right\}  \right)  \subseteq\mathcal{P}_{m}\left(
S\right)  $. In other words, $V$ is an $m$-element subset of $S$ (since
$\mathcal{P}_{m}\left(  S\right)  $ is the set of all $m$-element subsets of
$S$). Thus, $V$ is an $m$-element set. In other words, $\left\vert
V\right\vert =m$. Also, $V$ is a subset of $S$. In other words, $V\subseteq
S$.
\par
Notice that $W=\underbrace{V}_{\subseteq S}\setminus\left\{  s\right\}
\subseteq S\setminus\left\{  s\right\}  $. In other words, $W$ is a subset of
$S\setminus\left\{  s\right\}  $.
\par
Now, we shall prove that $s\in V$. Indeed, assume the contrary. Thus, $s\notin
V$. Hence, $V\setminus\left\{  s\right\}  =V$, so that $V=V\setminus\left\{
s\right\}  \subseteq S\setminus\left\{  s\right\}  $. Hence, $V$ is a subset
of $S\setminus\left\{  s\right\}  $. Thus, $V$ is an $m$-element subset of
$S\setminus\left\{  s\right\}  $. In other words, $V\in\mathcal{P}_{m}\left(
S\setminus\left\{  s\right\}  \right)  $ (since $\mathcal{P}_{m}\left(
S\setminus\left\{  s\right\}  \right)  $ is the set of all $m$-element subsets
of $S\setminus\left\{  s\right\}  $). This contradicts $V\notin\mathcal{P}%
_{m}\left(  S\setminus\left\{  s\right\}  \right)  $. This contradiction
proves that our assumption was wrong.
\par
Hence, $s\in V$ is proven. Thus, $\left\vert V\setminus\left\{  s\right\}
\right\vert =\underbrace{\left\vert V\right\vert }_{=m}-1=m-1$. In other
words, $V\setminus\left\{  s\right\}  $ is an $\left(  m-1\right)  $-element
set. In other words, $W$ is an $\left(  m-1\right)  $-element set (since
$W=V\setminus\left\{  s\right\}  $). Thus, $W$ is an $\left(  m-1\right)
$-element subset of $S\setminus\left\{  s\right\}  $. In other words,
$W\in\mathcal{P}_{m-1}\left(  S\setminus\left\{  s\right\}  \right)  $ (since
$\mathcal{P}_{m-1}\left(  S\setminus\left\{  s\right\}  \right)  $ is the set
of all $\left(  m-1\right)  $-element subsets of $S\setminus\left\{
s\right\}  $). Hence, $V\setminus\left\{  s\right\}  =W\in\mathcal{P}%
_{m-1}\left(  S\setminus\left\{  s\right\}  \right)  $. Qed.}. Hence, we can
define a map%
\[
\beta:\mathcal{P}_{m}\left(  S\right)  \setminus\mathcal{P}_{m}\left(
S\setminus\left\{  s\right\}  \right)  \rightarrow\mathcal{P}_{m-1}\left(
S\setminus\left\{  s\right\}  \right)
\]
by%
\[
\left(  \beta\left(  V\right)  =V\setminus\left\{  s\right\}
\ \ \ \ \ \ \ \ \ \ \text{for every }V\in\mathcal{P}_{m}\left(  S\right)
\setminus\mathcal{P}_{m}\left(  S\setminus\left\{  s\right\}  \right)
\right)  .
\]
Consider this map $\beta$.

We have $\alpha\circ\beta=\operatorname*{id}$\ \ \ \ \footnote{\textit{Proof.}
Let $V\in\mathcal{P}_{m}\left(  S\right)  \setminus\mathcal{P}_{m}\left(
S\setminus\left\{  s\right\}  \right)  $. We shall show that $\left(
\alpha\circ\beta\right)  \left(  V\right)  =\operatorname*{id}\left(
V\right)  $.
\par
We have $V\notin\mathcal{P}_{m}\left(  S\setminus\left\{  s\right\}  \right)
$ (since $V\in\mathcal{P}_{m}\left(  S\right)  \setminus\mathcal{P}_{m}\left(
S\setminus\left\{  s\right\}  \right)  $).
\par
We have $V\in\mathcal{P}_{m}\left(  S\right)  \setminus\mathcal{P}_{m}\left(
S\setminus\left\{  s\right\}  \right)  \subseteq\mathcal{P}_{m}\left(
S\right)  $. In other words, $V$ is an $m$-element subset of $S$ (since
$\mathcal{P}_{m}\left(  S\right)  $ is the set of all $m$-element subsets of
$S$). Thus, $V$ is a subset of $S$. In other words, $V\subseteq S$.
\par
Now, we shall prove that $s\in V$. Indeed, assume the contrary. Thus, $s\notin
V$. Hence, $V\setminus\left\{  s\right\}  =V$, so that $V=\underbrace{V}%
_{\subseteq S}\setminus\left\{  s\right\}  \subseteq S\setminus\left\{
s\right\}  $. Hence, $V$ is a subset of $S\setminus\left\{  s\right\}  $.
Thus, $V$ is an $m$-element subset of $S\setminus\left\{  s\right\}  $. In
other words, $V\in\mathcal{P}_{m}\left(  S\setminus\left\{  s\right\}
\right)  $ (since $\mathcal{P}_{m}\left(  S\setminus\left\{  s\right\}
\right)  $ is the set of all $m$-element subsets of $S\setminus\left\{
s\right\}  $). This contradicts $V\notin\mathcal{P}_{m}\left(  S\setminus
\left\{  s\right\}  \right)  $. This contradiction proves that our assumption
was wrong.
\par
Hence, $s\in V$ is proven. Thus, $V\cup\left\{  s\right\}  =V$. Now,
\begin{align*}
\left(  \alpha\circ\beta\right)  \left(  V\right)   &  =\alpha\left(
\beta\left(  V\right)  \right)  =\underbrace{\beta\left(  V\right)
}_{\substack{=V\setminus\left\{  s\right\}  \\\text{(by the definition of
}\beta\text{)}}}\cup\left\{  s\right\}  \ \ \ \ \ \ \ \ \ \ \left(  \text{by
the definition of }\alpha\right) \\
&  =\left(  V\setminus\left\{  s\right\}  \right)  \cup\left\{  s\right\}
=V\cup\left\{  s\right\}  =V=\operatorname*{id}\left(  V\right)  .
\end{align*}
\par
Now, forget that we fixed $V$. We thus have shown that $\left(  \alpha
\circ\beta\right)  \left(  V\right)  =\operatorname*{id}\left(  V\right)  $
for every $V\in\mathcal{P}_{m}\left(  S\right)  \setminus\mathcal{P}%
_{m}\left(  S\setminus\left\{  s\right\}  \right)  $. In other words,
$\alpha\circ\beta=\operatorname*{id}$. Qed.} and $\beta\circ\alpha
=\operatorname*{id}$\ \ \ \ \footnote{\textit{Proof.} Let $U\in\mathcal{P}%
_{m-1}\left(  S\setminus\left\{  s\right\}  \right)  $. We shall prove that
$\left(  \beta\circ\alpha\right)  \left(  U\right)  =\operatorname*{id}\left(
U\right)  $.
\par
We have $U\in\mathcal{P}_{m-1}\left(  S\setminus\left\{  s\right\}  \right)
$. In other words, $U$ is an $\left(  m-1\right)  $-element subset of
$S\setminus\left\{  s\right\}  $ (since $\mathcal{P}_{m-1}\left(
S\setminus\left\{  s\right\}  \right)  $ is the set of all $\left(
m-1\right)  $-element subsets of $S\setminus\left\{  s\right\}  $). Thus, $U$
is a subset of $S\setminus\left\{  s\right\}  $; hence, $U\subseteq
S\setminus\left\{  s\right\}  $.
\par
We have $s\in\left\{  s\right\}  $ and thus $s\notin S\setminus\left\{
s\right\}  $. If we had $s\in U$, then we would have $s\in U\subseteq
S\setminus\left\{  s\right\}  $, which would contradict $s\notin
S\setminus\left\{  s\right\}  $. Thus, we cannot have $s\in U$. In other
words, we have $s\notin U$. Hence, $U\setminus\left\{  s\right\}  =U$.
\par
Now,%
\begin{align*}
\left(  \beta\circ\alpha\right)  \left(  U\right)   &  =\beta\left(
\alpha\left(  U\right)  \right)  =\underbrace{\alpha\left(  U\right)
}_{\substack{=U\cup\left\{  s\right\}  \\\text{(by the definition of }%
\alpha\text{)}}}\setminus\left\{  s\right\}  \ \ \ \ \ \ \ \ \ \ \left(
\text{by the definition of }\beta\right) \\
&  =\left(  U\cup\left\{  s\right\}  \right)  \setminus\left\{  s\right\}
=U\setminus\left\{  s\right\}  =U=\operatorname*{id}\left(  U\right)  .
\end{align*}
\par
Now, forget that we fixed $U$. Thus, we have shown that $\left(  \beta
\circ\alpha\right)  \left(  U\right)  =\operatorname*{id}\left(  U\right)  $
for each $U\in\mathcal{P}_{m-1}\left(  S\setminus\left\{  s\right\}  \right)
$. In other words, $\beta\circ\alpha=\operatorname*{id}$, qed.}. These two
equalities show that the maps $\alpha$ and $\beta$ are mutually inverse. Thus,
the map $\alpha$ is invertible, i.e., is a bijection.

But $\alpha$ is the map
\begin{align*}
\mathcal{P}_{m-1}\left(  S\setminus\left\{  s\right\}  \right)   &
\rightarrow\mathcal{P}_{m}\left(  S\right)  \setminus\mathcal{P}_{m}\left(
S\setminus\left\{  s\right\}  \right)  ,\\
U  &  \mapsto U\cup\left\{  s\right\}
\end{align*}
(because $\alpha$ is the map $\mathcal{P}_{m-1}\left(  S\setminus\left\{
s\right\}  \right)  \rightarrow\mathcal{P}_{m}\left(  S\right)  \setminus
\mathcal{P}_{m}\left(  S\setminus\left\{  s\right\}  \right)  $ defined by
(\ref{pf.prop.sol.prop.binom.subsets.Pm.lem.c.alpha()=})).

We know that the map $\alpha$ is well-defined and a bijection. In other words,
the map
\begin{align*}
\mathcal{P}_{m-1}\left(  S\setminus\left\{  s\right\}  \right)   &
\rightarrow\mathcal{P}_{m}\left(  S\right)  \setminus\mathcal{P}_{m}\left(
S\setminus\left\{  s\right\}  \right)  ,\\
U  &  \mapsto U\cup\left\{  s\right\}
\end{align*}
is well-defined and a bijection (since $\alpha$ is the map
\begin{align*}
\mathcal{P}_{m-1}\left(  S\setminus\left\{  s\right\}  \right)   &
\rightarrow\mathcal{P}_{m}\left(  S\right)  \setminus\mathcal{P}_{m}\left(
S\setminus\left\{  s\right\}  \right)  ,\\
U  &  \mapsto U\cup\left\{  s\right\}
\end{align*}
). This proves Proposition \ref{prop.sol.prop.binom.subsets.Pm.lem}
\textbf{(c)}.
\end{proof}

We can now prove Proposition \ref{prop.binom.subsets} in a slightly rewritten form:

\begin{proposition}
\label{prop.sol.prop.binom.subsets.binom}Let $m\in\mathbb{N}$. Let $S$ be a
finite set. Then,%
\[
\left\vert \mathcal{P}_{m}\left(  S\right)  \right\vert =\dbinom{\left\vert
S\right\vert }{m}.
\]

\end{proposition}

\begin{proof}
[Proof of Proposition \ref{prop.sol.prop.binom.subsets.binom}.]We shall prove
Proposition \ref{prop.sol.prop.binom.subsets.binom} by induction over
$\left\vert S\right\vert $ (using the fact that $\left\vert S\right\vert
\in\mathbb{N}$ for any finite set $S$):

\textit{Induction base:} We must show that Proposition
\ref{prop.sol.prop.binom.subsets.binom} holds for $\left\vert S\right\vert
=0$. So, let $m\in\mathbb{N}$, and let $S$ be a finite set satisfying
$\left\vert S\right\vert =0$. We must then prove Proposition
\ref{prop.sol.prop.binom.subsets.binom}.

We have $\left\vert S\right\vert =0$ and thus $S=\varnothing$. We now want to
prove that%
\begin{equation}
\left\vert \mathcal{P}_{m}\left(  S\right)  \right\vert =\dbinom{\left\vert
S\right\vert }{m}. \label{eq.prop.sol.prop.binom.subsets.binom.ibas.goal}%
\end{equation}

We are in one of the following two cases:

\textit{Case 1:} We have $m=0$.

\textit{Case 2:} We have $m\neq0$.

Let us first consider Case 1. In this case, we have $m=0$. Hence,%
\[
\mathcal{P}_{m}\left(  S\right)  =\mathcal{P}_{0}\left(  S\right)  =\left\{
\varnothing\right\}  \ \ \ \ \ \ \ \ \ \ \left(  \text{by Proposition
\ref{prop.sol.prop.binom.subsets.Pm.lem} \textbf{(a)}}\right)
\]
and thus $\left\vert \underbrace{\mathcal{P}_{m}\left(  S\right)  }_{=\left\{
\varnothing\right\}  }\right\vert =\left\vert \left\{  \varnothing\right\}
\right\vert =1$. Comparing this with%
\begin{align*}
\dbinom{\left\vert S\right\vert }{m}  &  =\dbinom{\left\vert S\right\vert }%
{0}\ \ \ \ \ \ \ \ \ \ \left(  \text{since }m=0\right) \\
&  =1\ \ \ \ \ \ \ \ \ \ \left(  \text{by (\ref{eq.binom.00}), applied to
}\left\vert S\right\vert \text{ instead of }m\right)  ,
\end{align*}
we obtain $\left\vert \mathcal{P}_{m}\left(  S\right)  \right\vert
=\dbinom{\left\vert S\right\vert }{m}$. Thus,
(\ref{eq.prop.sol.prop.binom.subsets.binom.ibas.goal}) is proven in Case 1.

Let us now consider Case 2. In this case, we have $m\neq0$. Combined with
$m\in\mathbb{N}$, this yields $m\in\mathbb{N}\setminus\left\{  0\right\}
=\left\{  1,2,3,\ldots\right\}  $. In other words, $m$ is a positive integer.
Thus, $m>0$, so that $0<m$. Hence, (\ref{eq.binom.0}) (applied to $0$ and $m$
instead of $m$ and $n$) yields $\dbinom{0}{m}=0$. But Proposition
\ref{prop.sol.prop.binom.subsets.Pm.lem} \textbf{(b)} yields $\mathcal{P}%
_{m}\left(  \varnothing\right)  =\varnothing$. Now,
\[
\left\vert \mathcal{P}_{m}\left(  \underbrace{S}_{=\varnothing}\right)
\right\vert =\left\vert \underbrace{\mathcal{P}_{m}\left(  \varnothing\right)
}_{=\varnothing}\right\vert =\left\vert \varnothing\right\vert =0.
\]
Comparing this with%
\begin{align*}
\dbinom{\left\vert S\right\vert }{m}  &  =\dbinom{0}{m}%
\ \ \ \ \ \ \ \ \ \ \left(  \text{since }\left\vert S\right\vert =0\right) \\
&  =0,
\end{align*}
we obtain $\left\vert \mathcal{P}_{m}\left(  S\right)  \right\vert
=\dbinom{\left\vert S\right\vert }{m}$. Thus,
(\ref{eq.prop.sol.prop.binom.subsets.binom.ibas.goal}) is proven in Case 2.

We now have proven (\ref{eq.prop.sol.prop.binom.subsets.binom.ibas.goal}) in
each of the two Cases 1 and 2. Since these two Cases cover all possibilities,
this yields that (\ref{eq.prop.sol.prop.binom.subsets.binom.ibas.goal}) always
holds. In other words, we have $\left\vert \mathcal{P}_{m}\left(  S\right)
\right\vert =\dbinom{\left\vert S\right\vert }{m}$. In other words,
Proposition \ref{prop.sol.prop.binom.subsets.binom} holds. Thus, Proposition
\ref{prop.sol.prop.binom.subsets.binom} is proven for $\left\vert S\right\vert
=0$. This completes the induction base.

\textit{Induction step:} Let $N$ be a positive integer. Assume that
Proposition \ref{prop.sol.prop.binom.subsets.binom} holds for $\left\vert
S\right\vert =N-1$. We must then show that Proposition
\ref{prop.sol.prop.binom.subsets.binom} holds for $\left\vert S\right\vert =N$.

We have assumed that Proposition \ref{prop.sol.prop.binom.subsets.binom} holds
for $\left\vert S\right\vert =N-1$. In other words, every $m\in\mathbb{N}$ and
every finite set $S$ satisfying $\left\vert S\right\vert =N-1$ satisfy%
\begin{equation}
\left\vert \mathcal{P}_{m}\left(  S\right)  \right\vert =\dbinom{\left\vert
S\right\vert }{m}. \label{eq.prop.sol.prop.binom.subsets.binom.istep.hyp}%
\end{equation}

Let now $m\in\mathbb{N}$, and let $S$ be a finite set satisfying $\left\vert
S\right\vert =N$. We shall show that%
\begin{equation}
\left\vert \mathcal{P}_{m}\left(  S\right)  \right\vert =\dbinom{\left\vert
S\right\vert }{m}. \label{eq.prop.sol.prop.binom.subsets.binom.istep.goal}%
\end{equation}

[\textit{Proof of (\ref{eq.prop.sol.prop.binom.subsets.binom.istep.goal}):} We
have $\left\vert S\right\vert =N>0$ (since $N$ is positive). Hence, the set
$S$ is nonempty. In other words, the set $S$ has at least one element $s$. Fix
such an $s$.

We are in one of the following two cases:

\textit{Case 1:} We have $m=0$.

\textit{Case 2:} We have $m\neq0$.

Let us first consider Case 1. In this case, we have $m=0$. Hence,%
\[
\mathcal{P}_{m}\left(  S\right)  =\mathcal{P}_{0}\left(  S\right)  =\left\{
\varnothing\right\}  \ \ \ \ \ \ \ \ \ \ \left(  \text{by Proposition
\ref{prop.sol.prop.binom.subsets.Pm.lem} \textbf{(a)}}\right)
\]
and thus $\left\vert \underbrace{\mathcal{P}_{m}\left(  S\right)  }_{=\left\{
\varnothing\right\}  }\right\vert =\left\vert \left\{  \varnothing\right\}
\right\vert =1$. Comparing this with%
\begin{align*}
\dbinom{\left\vert S\right\vert }{m}  &  =\dbinom{\left\vert S\right\vert }%
{0}\ \ \ \ \ \ \ \ \ \ \left(  \text{since }m=0\right) \\
&  =1\ \ \ \ \ \ \ \ \ \ \left(  \text{by (\ref{eq.binom.00}), applied to
}\left\vert S\right\vert \text{ instead of }m\right)  ,
\end{align*}
we obtain $\left\vert \mathcal{P}_{m}\left(  S\right)  \right\vert
=\dbinom{\left\vert S\right\vert }{m}$. Thus,
(\ref{eq.prop.sol.prop.binom.subsets.binom.istep.goal}) is proven in Case 1.

Let us now consider Case 2. In this case, we have $m\neq0$. Combined with
$m\in\mathbb{N}$, this yields $m\in\mathbb{N}\setminus\left\{  0\right\}
=\left\{  1,2,3,\ldots\right\}  $. In other words, $m$ is a positive integer.
Hence, Proposition \ref{prop.sol.prop.binom.subsets.Pm.lem} \textbf{(c)}
yields that $\mathcal{P}_{m}\left(  S\setminus\left\{  s\right\}  \right)
\subseteq\mathcal{P}_{m}\left(  S\right)  $, and that the map%
\begin{align*}
\mathcal{P}_{m-1}\left(  S\setminus\left\{  s\right\}  \right)   &
\rightarrow\mathcal{P}_{m}\left(  S\right)  \setminus\mathcal{P}_{m}\left(
S\setminus\left\{  s\right\}  \right)  ,\\
U  &  \mapsto U\cup\left\{  s\right\}
\end{align*}
is well-defined and a bijection.

So there exists a bijection from $\mathcal{P}_{m-1}\left(  S\setminus\left\{
s\right\}  \right)  $ to $\mathcal{P}_{m}\left(  S\right)  \setminus
\mathcal{P}_{m}\left(  S\setminus\left\{  s\right\}  \right)  $ (namely, the
map
\begin{align*}
\mathcal{P}_{m-1}\left(  S\setminus\left\{  s\right\}  \right)   &
\rightarrow\mathcal{P}_{m}\left(  S\right)  \setminus\mathcal{P}_{m}\left(
S\setminus\left\{  s\right\}  \right)  ,\\
U  &  \mapsto U\cup\left\{  s\right\}
\end{align*}
). Hence, $\left\vert \mathcal{P}_{m}\left(  S\right)  \setminus
\mathcal{P}_{m}\left(  S\setminus\left\{  s\right\}  \right)  \right\vert
=\left\vert \mathcal{P}_{m-1}\left(  S\setminus\left\{  s\right\}  \right)
\right\vert $.

But $s\in S$ and thus $\left\vert S\setminus\left\{  s\right\}  \right\vert
=\underbrace{\left\vert S\right\vert }_{=N}-1=N-1$. Hence,
(\ref{eq.prop.sol.prop.binom.subsets.binom.istep.hyp}) (applied to
$S\setminus\left\{  s\right\}  $ instead of $S$) yields
\[
\left\vert \mathcal{P}_{m}\left(  S\setminus\left\{  s\right\}  \right)
\right\vert =\dbinom{\left\vert S\setminus\left\{  s\right\}  \right\vert }%
{m}=\dbinom{N-1}{m}\ \ \ \ \ \ \ \ \ \ \left(  \text{since }\left\vert
S\setminus\left\{  s\right\}  \right\vert =N-1\right)  .
\]
Also, $m$ is a positive integer; thus, $m-1\in\mathbb{N}$. Hence,
(\ref{eq.prop.sol.prop.binom.subsets.binom.istep.hyp}) (applied to $m-1$ and
$S\setminus\left\{  s\right\}  $ instead of $m$ and $S$) yields
\[
\left\vert \mathcal{P}_{m-1}\left(  S\setminus\left\{  s\right\}  \right)
\right\vert =\dbinom{\left\vert S\setminus\left\{  s\right\}  \right\vert
}{m-1}=\dbinom{N-1}{m-1}\ \ \ \ \ \ \ \ \ \ \left(  \text{since }\left\vert
S\setminus\left\{  s\right\}  \right\vert =N-1\right)  .
\]

Now, recall that $\mathcal{P}_{m}\left(  S\setminus\left\{  s\right\}
\right)  \subseteq\mathcal{P}_{m}\left(  S\right)  $. Thus,%
\begin{align*}
\left\vert \mathcal{P}_{m}\left(  S\right)  \right\vert  &
=\underbrace{\left\vert \mathcal{P}_{m}\left(  S\setminus\left\{  s\right\}
\right)  \right\vert }_{=\dbinom{N-1}{m}}+\underbrace{\left\vert
\mathcal{P}_{m}\left(  S\right)  \setminus\mathcal{P}_{m}\left(
S\setminus\left\{  s\right\}  \right)  \right\vert }_{\substack{=\left\vert
\mathcal{P}_{m-1}\left(  S\setminus\left\{  s\right\}  \right)  \right\vert
\\=\dbinom{N-1}{m-1}}}\\
&  =\dbinom{N-1}{m}+\dbinom{N-1}{m-1}=\dbinom{N-1}{m-1}+\dbinom{N-1}{m}.
\end{align*}
Comparing this with%
\begin{align*}
\dbinom{\left\vert S\right\vert }{m}  &  =\dbinom{N}{m}%
\ \ \ \ \ \ \ \ \ \ \left(  \text{since }\left\vert S\right\vert =N\right) \\
&  =\dbinom{N-1}{m-1}+\dbinom{N-1}{m}\ \ \ \ \ \ \ \ \ \ \left(
\begin{array}
[c]{c}%
\text{by (\ref{eq.binom.rec.m}), applied to }N\text{ and }m\\
\text{instead of }m\text{ and }n
\end{array}
\right)  ,
\end{align*}
we obtain $\left\vert \mathcal{P}_{m}\left(  S\right)  \right\vert
=\dbinom{\left\vert S\right\vert }{m}$. Thus,
(\ref{eq.prop.sol.prop.binom.subsets.binom.istep.goal}) is proven in Case 2.

We now have proven (\ref{eq.prop.sol.prop.binom.subsets.binom.istep.goal}) in
each of the two Cases 1 and 2. Since these two Cases cover all possibilities,
this yields that (\ref{eq.prop.sol.prop.binom.subsets.binom.istep.goal})
always holds.]

So we have proven (\ref{eq.prop.sol.prop.binom.subsets.binom.istep.goal}). In
other words, we have $\left\vert \mathcal{P}_{m}\left(  S\right)  \right\vert
=\dbinom{\left\vert S\right\vert }{m}$.

Now, forget that we fixed $m$ and $S$. We thus have shown that if
$m\in\mathbb{N}$, and if $S$ is a finite set satisfying $\left\vert
S\right\vert =N$, then $\left\vert \mathcal{P}_{m}\left(  S\right)
\right\vert =\dbinom{\left\vert S\right\vert }{m}$. In other words,
Proposition \ref{prop.sol.prop.binom.subsets.binom} is proven for $\left\vert
S\right\vert =N$. This completes the induction step. The induction proof of
Proposition \ref{prop.sol.prop.binom.subsets.binom} is thus complete.
\end{proof}

\begin{proof}
[Proof of Proposition \ref{prop.binom.subsets}.]Let $m\in\mathbb{N}$ and
$n\in\mathbb{N}$. Let $S$ be an $m$-element set. Thus, $S$ is a finite set.
Moreover, $\left\vert S\right\vert =m$ (since $S$ is an $m$-element set).

But $\mathcal{P}_{n}\left(  S\right)  $ is the set of all $n$-element subsets
of $S$ (by the definition of $\mathcal{P}_{n}\left(  S\right)  $). Hence,%
\[
\mathcal{P}_{n}\left(  S\right)  =\left(  \text{the set of all }%
n\text{-element subsets of }S\right)  .
\]
Thus,%
\begin{align*}
&  \left\vert \underbrace{\mathcal{P}_{n}\left(  S\right)  }_{=\left(
\text{the set of all }n\text{-element subsets of }S\right)  }\right\vert \\
&  =\left\vert \left(  \text{the set of all }n\text{-element subsets of
}S\right)  \right\vert \\
&  =\left(  \text{the number of all }n\text{-element subsets of }S\right)  .
\end{align*}
Comparing this with%
\begin{align*}
\left\vert \mathcal{P}_{n}\left(  S\right)  \right\vert  &  =\dbinom
{\left\vert S\right\vert }{n}\ \ \ \ \ \ \ \ \ \ \left(  \text{by Proposition
\ref{prop.sol.prop.binom.subsets.binom}, applied to }n\text{ instead of
}m\right) \\
&  =\dbinom{m}{n}\ \ \ \ \ \ \ \ \ \ \left(  \text{since }\left\vert
S\right\vert =m\right)  ,
\end{align*}
we obtain%
\[
\dbinom{m}{n}=\left(  \text{the number of all }n\text{-element subsets of
}S\right)  .
\]
In other words, $\dbinom{m}{n}$ is the number of all $n$-element subsets of
$S$. This proves Proposition \ref{prop.binom.subsets}.
\end{proof}
\end{verlong}

\subsection{Solution to Exercise \ref{exe.bin.-1and-2}}

Exercise \ref{exe.bin.-1and-2} \textbf{(a)} is precisely the statement of
Corollary \ref{cor.binom.-1} (with $n$ renamed as $k$). Let me still give two
proofs for it.

\begin{proof}
[Solution to Exercise \ref{exe.bin.-1and-2}.]\textbf{(a)} \textit{First
solution to Exercise \ref{exe.bin.-1and-2} \textbf{(a)}:} Let $k\in\mathbb{N}%
$. Then, Proposition \ref{prop.binom.upper-neg} (applied to $m=-1$ and $n=k$)
yields%
\begin{align*}
\dbinom{-1}{k}  &  =\left(  -1\right)  ^{k}\dbinom{k-\left(  -1\right)  -1}%
{k}\\
&  =\left(  -1\right)  ^{k}\underbrace{\dbinom{k}{k}}_{\substack{=1\\\text{(by
Proposition \ref{prop.binom.mm}}\\\text{(applied to }m=k\text{))}%
}}\ \ \ \ \ \ \ \ \ \ \left(  \text{since }k-\left(  -1\right)  -1=k\right) \\
&  =\left(  -1\right)  ^{k}.
\end{align*}
This solves Exercise \ref{exe.bin.-1and-2} \textbf{(a)}.

\textit{Second solution to Exercise \ref{exe.bin.-1and-2} \textbf{(a)}:} Let
$k\in\mathbb{N}$. Then, the definition of $k!$ yields $k!=1\cdot2\cdot
\cdots\cdot k$. But (\ref{eq.binom.mn}) (applied to $-1$ and $k$ instead of
$m$ and $n$) yields%
\begin{align*}
\dbinom{-1}{k}  &  =\dfrac{\left(  -1\right)  \left(  -2\right)  \cdots\left(
-1-k+1\right)  }{k!}\\
&  =\dfrac{\left(  -1\right)  \left(  -2\right)  \cdots\left(  -k\right)
}{k!}\ \ \ \ \ \ \ \ \ \ \left(  \text{since }-1-k+1=-k\right) \\
&  =\dfrac{1}{k!}\cdot\underbrace{\left(  -1\right)  \left(  -2\right)
\cdots\left(  -k\right)  }_{=\left(  -1\right)  ^{k}\cdot\left(  1\cdot
2\cdot\cdots\cdot k\right)  }=\dfrac{1}{k!}\cdot\left(  -1\right)  ^{k}%
\cdot\underbrace{\left(  1\cdot2\cdot\cdots\cdot k\right)  }_{=k!}\\
&  =\dfrac{1}{k!}\cdot\left(  -1\right)  ^{k}\cdot k!=\left(  -1\right)  ^{k}.
\end{align*}
This solves Exercise \ref{exe.bin.-1and-2} \textbf{(a)}.

\textbf{(b)} \textit{First solution to Exercise \ref{exe.bin.-1and-2}
\textbf{(b)}:} Let $k\in\mathbb{N}$. Then, $k+1\in\mathbb{N}$ and $k+1\geq k$.
Hence, Proposition \ref{prop.binom.symm} (applied to $m=k+1$ and $n=k$) yields%
\begin{align*}
\dbinom{k+1}{k}  &  =\dbinom{k+1}{\left(  k+1\right)  -k}=\dbinom{k+1}%
{1}\ \ \ \ \ \ \ \ \ \ \left(  \text{since }\left(  k+1\right)  -k=1\right) \\
&  =k+1\ \ \ \ \ \ \ \ \ \ \left(  \text{by Proposition \ref{prop.binom.00}
\textbf{(b)} (applied to }m=k+1\text{)}\right)  .
\end{align*}

But Proposition \ref{prop.binom.upper-neg} (applied to $m=-2$ and $n=k$)
yields%
\begin{align*}
\dbinom{-2}{k}  &  =\left(  -1\right)  ^{k}\dbinom{k-\left(  -2\right)  -1}%
{k}\\
&  =\left(  -1\right)  ^{k}\underbrace{\dbinom{k+1}{k}}_{=k+1}%
\ \ \ \ \ \ \ \ \ \ \left(  \text{since }k-\left(  -2\right)  -1=k+1\right) \\
&  =\left(  -1\right)  ^{k}\left(  k+1\right)  .
\end{align*}
This solves Exercise \ref{exe.bin.-1and-2} \textbf{(b)}.

\textit{Second solution to Exercise \ref{exe.bin.-1and-2} \textbf{(b)}:} Let
$k\in\mathbb{N}$. Then, the definition of $\left(  k+1\right)  !$ yields%
\[
\left(  k+1\right)  !=1\cdot2\cdot\cdots\cdot\left(  k+1\right)
=1\cdot\left(  2\cdot3\cdot\cdots\cdot\left(  k+1\right)  \right)
=2\cdot3\cdot\cdots\cdot\left(  k+1\right)  .
\]
On the other hand, applying (\ref{eq.n!.rec}) to $n=k+1$, we find%
\[
\left(  k+1\right)  !=\left(  k+1\right)  \cdot\left(  \underbrace{\left(
k+1\right)  -1}_{=k}\right)  !=\left(  k+1\right)  \cdot k!.
\]
But (\ref{eq.binom.mn}) (applied to $-2$ and $k$ instead of $m$ and $n$)
yields%
\begin{align*}
\dbinom{-2}{k}  &  =\dfrac{\left(  -2\right)  \left(  -3\right)  \cdots\left(
-2-k+1\right)  }{k!}\\
&  =\dfrac{\left(  -2\right)  \left(  -3\right)  \cdots\left(  -\left(
k+1\right)  \right)  }{k!}\ \ \ \ \ \ \ \ \ \ \left(  \text{since
}-2-k+1=-\left(  k+1\right)  \right) \\
&  =\dfrac{1}{k!}\cdot\underbrace{\left(  -2\right)  \left(  -3\right)
\cdots\left(  -\left(  k+1\right)  \right)  }_{=\left(  -1\right)  ^{k}%
\cdot\left(  2\cdot3\cdot\cdots\cdot\left(  k+1\right)  \right)  }=\dfrac
{1}{k!}\cdot\left(  -1\right)  ^{k}\cdot\underbrace{\left(  2\cdot3\cdot
\cdots\cdot\left(  k+1\right)  \right)  }_{=\left(  k+1\right)  !}\\
&  =\dfrac{1}{k!}\cdot\left(  -1\right)  ^{k}\cdot\underbrace{\left(
k+1\right)  !}_{=\left(  k+1\right)  \cdot k!}=\dfrac{1}{k!}\cdot\left(
-1\right)  ^{k}\cdot\left(  k+1\right)  \cdot k!=\left(  -1\right)
^{k}\left(  k+1\right)  .
\end{align*}
This solves Exercise \ref{exe.bin.-1and-2} \textbf{(b)}.

\textbf{(c)} We shall give two solutions to Exercise \ref{exe.bin.-1and-2}
\textbf{(c)}: one by astutely transforming the left-hand side, and another by
straightforward induction.

\textit{First solution to Exercise \ref{exe.bin.-1and-2} \textbf{(c)}:} Let
$n\in\mathbb{N}$. Then, we can group the factors of the product $1!\cdot
2!\cdot\cdots\cdot\left(  2n\right)  !$ into pairs of successive factors. We
thus obtain%
\begin{align*}
&  1!\cdot2!\cdot\cdots\cdot\left(  2n\right)  !\\
&  =\left(  1!\cdot2!\right)  \cdot\left(  3!\cdot4!\right)  \cdot\cdots
\cdot\left(  \left(  2n-1\right)  !\cdot\left(  2n\right)  !\right)
=\prod_{i=1}^{n}\left(  \left(  2i-1\right)  !\cdot\underbrace{\left(
2i\right)  !}_{\substack{=2i\cdot\left(  2i-1\right)  !\\\text{(by
(\ref{eq.n!.rec}), applied}\\\text{to }n=2i\text{)}}}\right) \\
&  =\prod_{i=1}^{n}\underbrace{\left(  \left(  2i-1\right)  !\cdot
2i\cdot\left(  2i-1\right)  !\right)  }_{=2i\cdot\left(  2i-1\right)  !^{2}%
}=\prod_{i=1}^{n}\left(  2i\cdot\left(  2i-1\right)  !^{2}\right) \\
&  =\underbrace{\left(  \prod_{i=1}^{n}\left(  2i\right)  \right)  }%
_{=2^{n}\prod_{i=1}^{n}i}\cdot\underbrace{\left(  \prod_{i=1}^{n}\left(
\left(  2i-1\right)  !^{2}\right)  \right)  }_{=\left(  \prod_{i=1}^{n}\left(
\left(  2i-1\right)  !\right)  \right)  ^{2}}=2^{n}\underbrace{\left(
\prod_{i=1}^{n}i\right)  }_{=1\cdot2\cdot\cdots\cdot n=n!}\cdot\left(
\prod_{i=1}^{n}\left(  \left(  2i-1\right)  !\right)  \right)  ^{2}\\
&  =2^{n}n!\cdot\left(  \prod_{i=1}^{n}\left(  \left(  2i-1\right)  !\right)
\right)  ^{2}.
\end{align*}
Dividing both sides of this equality by $n!$, we find%
\[
\dfrac{1!\cdot2!\cdot\cdots\cdot\left(  2n\right)  !}{n!}=2^{n}\cdot\left(
\prod_{i=1}^{n}\left(  \left(  2i-1\right)  !\right)  \right)  ^{2}.
\]
This solves Exercise \ref{exe.bin.-1and-2} \textbf{(c)}.

\textit{Second solution to Exercise \ref{exe.bin.-1and-2} \textbf{(c)}:} We
shall solve Exercise \ref{exe.bin.-1and-2} \textbf{(c)} by induction on $n$:

\textit{Induction base:} Comparing%
\[
\dfrac{1!\cdot2!\cdot\cdots\cdot\left(  2\cdot0\right)  !}{0!}%
=\underbrace{\dfrac{1}{0!}}_{=1}\cdot\underbrace{\left(  1!\cdot2!\cdot
\cdots\cdot\left(  2\cdot0\right)  !\right)  }_{\substack{=1!\cdot
2!\cdot\cdots\cdot0!\\=\left(  \text{empty product}\right)  =1}}=1
\]
with%
\[
\underbrace{2^{0}}_{=1}\cdot\left(  \underbrace{\prod_{i=1}^{0}\left(  \left(
2i-1\right)  !\right)  }_{=\left(  \text{empty product}\right)  =1}\right)
^{2}=1\cdot1^{2}=1,
\]
we obtain $\dfrac{1!\cdot2!\cdot\cdots\cdot\left(  2\cdot0\right)  !}%
{0!}=2^{0}\cdot\left(  \prod_{i=1}^{0}\left(  \left(  2i-1\right)  !\right)
\right)  ^{2}$. In other words, Exercise \ref{exe.bin.-1and-2} \textbf{(c)}
holds for $n=0$. This completes the induction base.

\textit{Induction step:} Let $m\in\mathbb{N}$. Assume that Exercise
\ref{exe.bin.-1and-2} \textbf{(c)} holds for $n=m$. We must prove that
Exercise \ref{exe.bin.-1and-2} \textbf{(c)} holds for $n=m+1$.

Clearly, $m+1$ is a positive integer (since $m+1\geq1>0$). Thus, applying
(\ref{eq.n!.rec}) to $n=m+1$, we find
\[
\left(  m+1\right)  !=\left(  m+1\right)  \cdot\left(  \underbrace{\left(
m+1\right)  -1}_{=m}\right)  !=\left(  m+1\right)  \cdot m!.
\]

On the other hand, $2m+2\geq2>0$, so that $2m+2$ is a positive integer. Hence,
applying (\ref{eq.n!.rec}) to $n=2m+2$, we find%
\[
\left(  2m+2\right)  !=\left(  2m+2\right)  \cdot\left(  \underbrace{\left(
2m+2\right)  -1}_{=2m+1}\right)  !=\left(  2m+2\right)  \cdot\left(
2m+1\right)  !.
\]

We have assumed that Exercise \ref{exe.bin.-1and-2} \textbf{(c)} holds for
$n=m$. In other words, we have%
\[
\dfrac{1!\cdot2!\cdot\cdots\cdot\left(  2m\right)  !}{m!}=2^{m}\cdot\left(
\prod_{i=1}^{m}\left(  \left(  2i-1\right)  !\right)  \right)  ^{2}.
\]
Now,%
\begin{align*}
&  \dfrac{1!\cdot2!\cdot\cdots\cdot\left(  2\left(  m+1\right)  \right)
!}{\left(  m+1\right)  !}\\
&  =\dfrac{1!\cdot2!\cdot\cdots\cdot\left(  2m+2\right)  !}{\left(
m+1\right)  \cdot m!}\\
&  \ \ \ \ \ \ \ \ \ \ \left(  \text{since }2\left(  m+1\right)  =2m+2\text{
and }\left(  m+1\right)  !=\left(  m+1\right)  \cdot m!\right) \\
&  =\dfrac{1}{\left(  m+1\right)  \cdot m!}\cdot\underbrace{\left(
1!\cdot2!\cdot\cdots\cdot\left(  2m+2\right)  !\right)  }_{=\left(
1!\cdot2!\cdot\cdots\cdot\left(  2m\right)  !\right)  \cdot\left(
2m+1\right)  !\cdot\left(  2m+2\right)  !}\\
&  =\dfrac{1}{\left(  m+1\right)  \cdot m!}\cdot\left(  1!\cdot2!\cdot
\cdots\cdot\left(  2m\right)  !\right)  \cdot\left(  2m+1\right)
!\cdot\underbrace{\left(  2m+2\right)  !}_{=\left(  2m+2\right)  \cdot\left(
2m+1\right)  !}\\
&  =\dfrac{1}{\left(  m+1\right)  \cdot m!}\cdot\left(  1!\cdot2!\cdot
\cdots\cdot\left(  2m\right)  !\right)  \cdot\left(  2m+1\right)
!\cdot\left(  2m+2\right)  \cdot\left(  2m+1\right)  !\\
&  =\underbrace{\dfrac{2m+2}{m+1}}_{=2}\cdot\underbrace{\dfrac{1!\cdot
2!\cdot\cdots\cdot\left(  2m\right)  !}{m!}}_{=2^{m}\cdot\left(  \prod
_{i=1}^{m}\left(  \left(  2i-1\right)  !\right)  \right)  ^{2}}\cdot\left(
\left(  2m+1\right)  !\right)  ^{2}\\
&  =2\cdot2^{m}\cdot\left(  \prod_{i=1}^{m}\left(  \left(  2i-1\right)
!\right)  \right)  ^{2}\cdot\left(  \left(  2m+1\right)  !\right)  ^{2}.
\end{align*}
Comparing this with%
\begin{align*}
&  \underbrace{2^{m+1}}_{=2\cdot2^{m}}\cdot\left(  \underbrace{\prod
_{i=1}^{m+1}\left(  \left(  2i-1\right)  !\right)  }_{=\left(  \prod_{i=1}%
^{m}\left(  \left(  2i-1\right)  !\right)  \right)  \cdot\left(  2\left(
m+1\right)  -1\right)  !}\right)  ^{2}\\
&  =2\cdot2^{m}\cdot\left(  \left(  \prod_{i=1}^{m}\left(  \left(
2i-1\right)  !\right)  \right)  \cdot\left(  \underbrace{2\left(  m+1\right)
-1}_{=2m+1}\right)  !\right)  ^{2}\\
&  =2\cdot2^{m}\cdot\left(  \left(  \prod_{i=1}^{m}\left(  \left(
2i-1\right)  !\right)  \right)  \cdot\left(  2m+1\right)  !\right)  ^{2}\\
&  =2\cdot2^{m}\cdot\left(  \prod_{i=1}^{m}\left(  \left(  2i-1\right)
!\right)  \right)  ^{2}\cdot\left(  \left(  2m+1\right)  !\right)  ^{2},
\end{align*}
we obtain%
\[
\dfrac{1!\cdot2!\cdot\cdots\cdot\left(  2\left(  m+1\right)  \right)
!}{\left(  m+1\right)  !}=2^{m+1}\cdot\left(  \prod_{i=1}^{m+1}\left(  \left(
2i-1\right)  !\right)  \right)  ^{2}.
\]
In other words, Exercise \ref{exe.bin.-1and-2} \textbf{(c)} holds for $n=m+1$.
This completes the induction step. Thus, Exercise \ref{exe.bin.-1and-2}
\textbf{(c)} is solved.
\end{proof}

\subsection{Solution to Exercise \ref{exe.prop.binom.binomial}}

In order to solve Exercise \ref{exe.prop.binom.binomial}, we need to prove
Proposition \ref{prop.binom.binomial}.

\begin{proof}
[Proof of Proposition \ref{prop.binom.binomial}.]We shall prove Proposition
\ref{prop.binom.binomial} by induction over $n$:

\textit{Induction base:} We have%
\[
\sum_{k=0}^{0}\dbinom{0}{k}x^{k}y^{0-k}=\underbrace{\dbinom{0}{0}%
}_{\substack{=1\\\text{(by (\ref{eq.binom.00}), applied to }m=0\text{)}%
}}\underbrace{x^{0}}_{=1}\underbrace{y^{0-0}}_{=y^{0}=1}=1.
\]
Comparing this with $\left(  x+y\right)  ^{0}=1$, we obtain $\left(
x+y\right)  ^{0}=\sum_{k=0}^{0}\dbinom{0}{k}x^{k}y^{0-k}$. In other words,
Proposition \ref{prop.binom.binomial} holds for $n=0$. This completes the
induction base.

\textit{Induction step:} Let $N$ be a positive integer. Assume that
Proposition \ref{prop.binom.binomial} holds for $n=N-1$. We must now prove
that Proposition \ref{prop.binom.binomial} holds for $n=N$.

Notice that $N\geq1$ (since $N$ is a positive integer), so that $N-1\geq0$.

We have assumed that Proposition \ref{prop.binom.binomial} holds for $n=N-1$.
In other words, we have%
\begin{equation}
\left(  x+y\right)  ^{N-1}=\sum_{k=0}^{N-1}\dbinom{N-1}{k}x^{k}y^{\left(
N-1\right)  -k}. \label{pf.prop.binom.binomial.ihyp}%
\end{equation}

Now,%
\begin{align}
&  \left(  x+y\right)  ^{N}\nonumber\\
&  =\left(  x+y\right)  \underbrace{\left(  x+y\right)  ^{N-1}}%
_{\substack{=\sum_{k=0}^{N-1}\dbinom{N-1}{k}x^{k}y^{\left(  N-1\right)
-k}\\\text{(by (\ref{pf.prop.binom.binomial.ihyp}))}}}\nonumber\\
&  =\left(  x+y\right)  \left(  \sum_{k=0}^{N-1}\dbinom{N-1}{k}x^{k}y^{\left(
N-1\right)  -k}\right) \nonumber\\
&  =\underbrace{x\sum_{k=0}^{N-1}\dbinom{N-1}{k}x^{k}y^{\left(  N-1\right)
-k}}_{=\sum_{k=0}^{N-1}\dbinom{N-1}{k}xx^{k}y^{\left(  N-1\right)  -k}%
}+\underbrace{y\sum_{k=0}^{N-1}\dbinom{N-1}{k}x^{k}y^{\left(  N-1\right)  -k}%
}_{=\sum_{k=0}^{N-1}\dbinom{N-1}{k}x^{k}yy^{\left(  N-1\right)  -k}%
}\nonumber\\
&  =\sum_{k=0}^{N-1}\underbrace{\dbinom{N-1}{k}}_{\substack{=\dbinom
{N-1}{\left(  k+1\right)  -1}\\\text{(since }k=\left(  k+1\right)  -1\text{)}%
}}\underbrace{xx^{k}}_{=x^{k+1}}\underbrace{y^{\left(  N-1\right)  -k}%
}_{\substack{=y^{N-\left(  k+1\right)  }\\\text{(since }\left(  N-1\right)
-k=N-\left(  k+1\right)  \text{)}}}\nonumber\\
&  \ \ \ \ \ \ \ \ \ \ +\sum_{k=0}^{N-1}\dbinom{N-1}{k}x^{k}%
\underbrace{yy^{\left(  N-1\right)  -k}}_{\substack{=y^{\left(  \left(
N-1\right)  -k\right)  +1}=y^{N-k}\\\text{(since }\left(  \left(  N-1\right)
-k\right)  +1=N-k\text{)}}}\nonumber\\
&  =\underbrace{\sum_{k=0}^{N-1}\dbinom{N-1}{\left(  k+1\right)  -1}%
x^{k+1}y^{N-\left(  k+1\right)  }}_{\substack{=\sum_{k=1}^{N}\dbinom{N-1}%
{k-1}x^{k}y^{N-k}\\\text{(here, we have substituted }k\text{ for }k+1\text{ in
the sum)}}}+\sum_{k=0}^{N-1}\dbinom{N-1}{k}x^{k}y^{N-k}\nonumber\\
&  =\sum_{k=1}^{N}\dbinom{N-1}{k-1}x^{k}y^{N-k}+\sum_{k=0}^{N-1}\dbinom
{N-1}{k}x^{k}y^{N-k}. \label{pf.prop.binom.binomial.1}%
\end{align}
But recall that $N-1\geq0$, so that $N-1\in\mathbb{N}$. Clearly, $N-1<N$.
Hence, (\ref{eq.binom.0}) (applied to $m=N-1$ and $n=N$) yields $\dbinom
{N-1}{N}=0$.

On the other hand, $N$ is a positive integer. Hence, $N\in\left\{
1,2,\ldots,N\right\}  $. Thus, we can split off the addend for $k=N$ from the
sum $\sum_{k=1}^{N}\dbinom{N-1}{k}x^{k}y^{N-k}$. We thus obtain%
\begin{align}
\sum_{k=1}^{N}\dbinom{N-1}{k}x^{k}y^{N-k}  &  =\sum_{k=1}^{N-1}\dbinom{N-1}%
{k}x^{k}y^{N-k}+\underbrace{\dbinom{N-1}{N}}_{=0}x^{N}y^{N-N}\nonumber\\
&  =\sum_{k=1}^{N-1}\dbinom{N-1}{k}x^{k}y^{N-k}+\underbrace{0x^{N}y^{N-N}%
}_{=0}\nonumber\\
&  =\sum_{k=1}^{N-1}\dbinom{N-1}{k}x^{k}y^{N-k}.
\label{pf.prop.binom.binomial.4b1}%
\end{align}

Furthermore, $N-1\geq0$, so that $0\in\left\{  0,1,\ldots,N-1\right\}  $.
Hence, we can split off the addend for $k=0$ from the sum $\sum_{k=0}%
^{N-1}\dbinom{N-1}{k}x^{k}y^{N-k}$. We thus obtain%
\begin{align*}
&  \sum_{k=0}^{N-1}\dbinom{N-1}{k}x^{k}y^{N-k}\\
&  =\underbrace{\dbinom{N-1}{0}}_{\substack{=1\\\text{(by (\ref{eq.binom.00}),
applied to }m=N-1\text{)}}}x^{0}y^{N-0}+\sum_{k=1}^{N-1}\dbinom{N-1}{k}%
x^{k}y^{N-k}\\
&  =x^{0}y^{N-0}+\sum_{k=1}^{N-1}\dbinom{N-1}{k}x^{k}y^{N-k}.
\end{align*}
Comparing this with%
\begin{align*}
&  \underbrace{\dbinom{N}{0}}_{\substack{=1\\\text{(by (\ref{eq.binom.00}),
applied to }m=N\text{)}}}x^{0}y^{N-0}+\underbrace{\sum_{k=1}^{N}\dbinom
{N-1}{k}x^{k}y^{N-k}}_{\substack{=\sum_{k=1}^{N-1}\dbinom{N-1}{k}x^{k}%
y^{N-k}\\\text{(by (\ref{pf.prop.binom.binomial.4b1}))}}}\\
&  =x^{0}y^{N-0}+\sum_{k=1}^{N-1}\dbinom{N-1}{k}x^{k}y^{N-k},
\end{align*}
we obtain%
\begin{equation}
\sum_{k=0}^{N-1}\dbinom{N-1}{k}x^{k}y^{N-k}=\dbinom{N}{0}x^{0}y^{N-0}%
+\sum_{k=1}^{N}\dbinom{N-1}{k}x^{k}y^{N-k}. \label{pf.prop.binom.binomial.4b}%
\end{equation}
Now, (\ref{pf.prop.binom.binomial.1}) becomes%
\begin{align*}
&  \left(  x+y\right)  ^{N}\\
&  =\sum_{k=1}^{N}\dbinom{N-1}{k-1}x^{k}y^{N-k}+\underbrace{\sum_{k=0}%
^{N-1}\dbinom{N-1}{k}x^{k}y^{N-k}}_{\substack{=\dbinom{N}{0}x^{0}y^{N-0}%
+\sum_{k=1}^{N}\dbinom{N-1}{k}x^{k}y^{N-k}\\\text{(by
(\ref{pf.prop.binom.binomial.4b}))}}}\\
&  =\sum_{k=1}^{N}\dbinom{N-1}{k-1}x^{k}y^{N-k}+\dbinom{N}{0}x^{0}y^{N-0}%
+\sum_{k=1}^{N}\dbinom{N-1}{k}x^{k}y^{N-k}\\
&  =\dbinom{N}{0}x^{0}y^{N-0}+\underbrace{\sum_{k=1}^{N}\dbinom{N-1}{k-1}%
x^{k}y^{N-k}+\sum_{k=1}^{N}\dbinom{N-1}{k}x^{k}y^{N-k}}_{=\sum_{k=1}%
^{N}\left(  \dbinom{N-1}{k-1}+\dbinom{N-1}{k}\right)  x^{k}y^{N-k}}\\
&  =\dbinom{N}{0}x^{0}y^{N-0}+\sum_{k=1}^{N}\left(  \dbinom{N-1}{k-1}%
+\dbinom{N-1}{k}\right)  x^{k}y^{N-k}.
\end{align*}
Comparing this with%
\begin{align*}
&  \sum_{k=0}^{N}\dbinom{N}{k}x^{k}y^{N-k}\\
&  =\dbinom{N}{0}x^{0}y^{N-0}+\sum_{k=1}^{N}\underbrace{\dbinom{N}{k}%
}_{\substack{=\dbinom{N-1}{k-1}+\dbinom{N-1}{k}\\\text{(by
(\ref{eq.binom.rec.m}), applied to }m=N\text{ and }n=k\text{)}}}x^{k}y^{N-k}\\
&  \ \ \ \ \ \ \ \ \ \ \left(  \text{here, we have split off the addend for
}k=0\text{ from the sum}\right) \\
&  =\dbinom{N}{0}x^{0}y^{N-0}+\sum_{k=1}^{N}\left(  \dbinom{N-1}{k-1}%
+\dbinom{N-1}{k}\right)  x^{k}y^{N-k},
\end{align*}
we obtain $\left(  x+y\right)  ^{N}=\sum_{k=0}^{N}\dbinom{N}{k}x^{k}y^{N-k}$.
In other words, Proposition \ref{prop.binom.binomial} holds for $n=N$. Hence,
Proposition \ref{prop.binom.binomial} is proven by induction.
\end{proof}

\subsection{Solution to Exercise \ref{exe.binom.scaryfrac}}

Exercise \ref{exe.binom.scaryfrac} may look scary, but it is a straightforward
exercise on induction (on $b$). To make our life a little bit easier, we shall
slightly relax the condition $a\leq b$ to $b\geq a-1$ (so that we can use the
case $b=a-1$ instead of $b=a$ as an induction base):

\begin{proposition}
\label{prop.binom.scaryfrac}Let $k$ be a positive integer. Let $a$ be a
positive integer such that $k\leq a$. Let $b\in\left\{  a-1,a,a+1,\ldots
\right\}  $. Then,%
\[
\dfrac{k-1}{k}\sum_{n=a}^{b}\dfrac{1}{\dbinom{n}{k}}=\dfrac{1}{\dbinom
{a-1}{k-1}}-\dfrac{1}{\dbinom{b}{k-1}}.
\]
(In particular, all fractions appearing in this equality are well-defined.)
\end{proposition}

Before we prove this proposition (which rather obviously encompasses the claim
of Exercise \ref{exe.binom.scaryfrac}), let us show a few lemmas:

\begin{lemma}
\label{lem.sol.binom.scaryfrac.absorp2}Let $m\in\mathbb{Q}$ and $n\in\left\{
1,2,3,\ldots\right\}  $. Then,%
\[
\dbinom{m}{n}=\dfrac{m-n+1}{n}\dbinom{m}{n-1}.
\]

\end{lemma}

\begin{proof}
[Proof of Lemma \ref{lem.sol.binom.scaryfrac.absorp2}.]We have $n\neq0$ (since
$n\in\left\{  1,2,3,\ldots\right\}  $). Hence, the fraction $\dfrac{m-n+1}{n}$
is well-defined.

We have $n\in\left\{  1,2,3,\ldots\right\}  $; in other words, $n$ is a
positive integer. Hence, $n!=n\cdot\left(  n-1\right)  !$.

We have $n\in\left\{  1,2,3,\ldots\right\}  \subseteq\mathbb{N}$. Thus, the
definition of $\dbinom{m}{n}$ yields%
\begin{align}
\dbinom{m}{n}  &  =\dfrac{m\left(  m-1\right)  \cdots\left(  m-n+1\right)
}{n!}\nonumber\\
&  =\dfrac{m\left(  m-1\right)  \cdots\left(  m-n+1\right)  }{n\cdot\left(
n-1\right)  !}\ \ \ \ \ \ \ \ \ \ \left(  \text{since }n!=n\cdot\left(
n-1\right)  !\right) \nonumber\\
&  =\dfrac{1}{n\cdot\left(  n-1\right)  !}\cdot\underbrace{\left(  m\left(
m-1\right)  \cdots\left(  m-n+1\right)  \right)  }_{\substack{=\left(
m\left(  m-1\right)  \cdots\left(  m-n+2\right)  \right)  \cdot\left(
m-n+1\right)  \\\text{(since }n\text{ is a positive integer)}}}\nonumber\\
&  =\dfrac{1}{n\cdot\left(  n-1\right)  !}\cdot\left(  m\left(  m-1\right)
\cdots\left(  m-n+2\right)  \right)  \cdot\left(  m-n+1\right) \nonumber\\
&  =\dfrac{m-n+1}{n}\cdot\dfrac{m\left(  m-1\right)  \cdots\left(
m-n+2\right)  }{\left(  n-1\right)  !}.
\label{pf.lem.sol.binom.scaryfrac.absorp2.1}%
\end{align}

Moreover, $n-1\in\mathbb{N}$ (since $n\in\left\{  1,2,3,\ldots\right\}  $), so
that the definition of $\dbinom{m}{n-1}$ yields%
\begin{align*}
\dbinom{m}{n-1}  &  =\dfrac{m\left(  m-1\right)  \cdots\left(  m-\left(
n-1\right)  +1\right)  }{\left(  n-1\right)  !}=\dfrac{m\left(  m-1\right)
\cdots\left(  m-n+2\right)  }{\left(  n-1\right)  !}\\
&  \ \ \ \ \ \ \ \ \ \ \left(  \text{since }m-\left(  n-1\right)
+1=m-n+2\right)  .
\end{align*}
Multiplying this equality by $\dfrac{m-n+1}{n}$, we obtain%
\[
\dfrac{m-n+1}{n}\dbinom{m}{n-1}=\dfrac{m-n+1}{n}\cdot\dfrac{m\left(
m-1\right)  \cdots\left(  m-n+2\right)  }{\left(  n-1\right)  !}.
\]
Comparing this with (\ref{pf.lem.sol.binom.scaryfrac.absorp2.1}), we obtain
$\dbinom{m}{n}=\dfrac{m-n+1}{n}\dbinom{m}{n-1}$. This proves Lemma
\ref{lem.sol.binom.scaryfrac.absorp2}.
\end{proof}

\begin{lemma}
\label{lem.sol.binom.scaryfrac.wd}Let $k$ be a positive integer. Let $a$ be a
positive integer such that $k\leq a$. Let $b\in\left\{  a-1,a,a+1,\ldots
\right\}  $.

\textbf{(a)} The fractions $\dfrac{k-1}{k}$, $\dfrac{1}{\dbinom{a-1}{k-1}}$
and $\dfrac{1}{\dbinom{b}{k-1}}$ are well-defined.

\textbf{(b)} For each $n\in\left\{  a,a+1,\ldots,b\right\}  $, the fraction
$\dfrac{1}{\dbinom{n}{k}}$ is well-defined.
\end{lemma}

\begin{proof}
[Proof of Lemma \ref{lem.sol.binom.scaryfrac.wd}.]We first observe that%
\begin{equation}
\dbinom{m}{n}\neq0 \label{pf.lem.sol.binom.scaryfrac.wd.neq0}%
\end{equation}
for any $m\in\mathbb{N}$ and $n\in\mathbb{N}$ satisfying $m\geq n$.

[\textit{Proof of (\ref{pf.lem.sol.binom.scaryfrac.wd.neq0}):} Let
$m\in\mathbb{N}$ and $n\in\mathbb{N}$ be such that $m\geq n$. Thus,
Proposition \ref{prop.binom.formula} yields $\dbinom{m}{n}=\dfrac
{m!}{n!\left(  m-n\right)  !}\neq0$ (since $m!\neq0$ (since $m!$ is a positive
integer)). This proves (\ref{pf.lem.sol.binom.scaryfrac.wd.neq0}).]

We have $k\neq0$ (since $k$ is a positive integer); thus, the fraction
$\dfrac{k-1}{k}$ is well-defined.

Also, $a-1\in\mathbb{N}$ (since $a$ is a positive integer) and $k-1\in
\mathbb{N}$ (since $k$ is a positive integer) and $a-1\geq k-1$ (since
$\underbrace{k}_{\leq a}-1\leq a-1$). Thus,
(\ref{pf.lem.sol.binom.scaryfrac.wd.neq0}) (applied to $a-1$ and $k-1$ instead
of $m$ and $n$) yields $\dbinom{a-1}{k-1}\neq0$. Hence, the fraction
$\dfrac{1}{\dbinom{a-1}{k-1}}$ is well-defined.

Also, $b\in\left\{  a-1,a,a+1,\ldots\right\}  $, so that $b\geq a-1\geq
k-1\geq0$ (since $k-1\in\mathbb{N}$). Hence, $b\in\mathbb{N}$ (since $b\geq0$
and $b\in\left\{  a-1,a,a+1,\ldots\right\}  \subseteq\mathbb{Z}$) and
$k-1\in\mathbb{N}$ and $b\geq k-1$. Thus,
(\ref{pf.lem.sol.binom.scaryfrac.wd.neq0}) (applied to $b$ and $k-1$ instead
of $m$ and $n$) yields $\dbinom{b}{k-1}\neq0$. Hence, the fraction $\dfrac
{1}{\dbinom{b}{k-1}}$ is well-defined.

We have now shown that the fractions $\dfrac{k-1}{k}$, $\dfrac{1}{\dbinom
{a-1}{k-1}}$ and $\dfrac{1}{\dbinom{b}{k-1}}$ are well-defined. This proves
Lemma \ref{lem.sol.binom.scaryfrac.wd} \textbf{(a)}.

\textbf{(b)} Let $n\in\left\{  a,a+1,\ldots,b\right\}  $. Thus, $n\geq
a>a-1\geq0$, so that $n\in\mathbb{N}$ (since $n\in\left\{  a,a+1,\ldots
,b\right\}  \subseteq\mathbb{Z}$). Also, $k\in\mathbb{N}$. Furthermore, $n\geq
a\geq k$ (since $k\leq a$). Hence, (\ref{pf.lem.sol.binom.scaryfrac.wd.neq0})
(applied to $n$ and $k$ instead of $m$ and $n$) yields $\dbinom{n}{k}\neq0$.
Hence, the fraction $\dfrac{1}{\dbinom{n}{k}}$ is well-defined. This proves
Lemma \ref{lem.sol.binom.scaryfrac.wd} \textbf{(b)}.
\end{proof}

\begin{proof}
[Proof of Proposition \ref{prop.binom.scaryfrac}.]All fractions appearing in
Proposition \ref{prop.binom.scaryfrac} are well-defined (because of Lemma
\ref{lem.sol.binom.scaryfrac.wd}). It thus remains to prove that%
\begin{equation}
\dfrac{k-1}{k}\sum_{n=a}^{b}\dfrac{1}{\dbinom{n}{k}}=\dfrac{1}{\dbinom
{a-1}{k-1}}-\dfrac{1}{\dbinom{b}{k-1}}. \label{pf.prop.binom.scaryfrac.goal}%
\end{equation}

Forget that we fixed $b$. We thus need to prove that
(\ref{pf.prop.binom.scaryfrac.goal}) holds for every $b\in\left\{
a-1,a,a+1,\ldots\right\}  $. We shall prove this by induction on $b$:

\textit{Induction base:} Comparing%
\[
\dfrac{k-1}{k}\underbrace{\sum_{n=a}^{a-1}\dfrac{1}{\dbinom{n}{k}}}_{=\left(
\text{empty sum}\right)  =0}=\dfrac{k-1}{k}\cdot0=0
\]
with $\dfrac{1}{\dbinom{a-1}{k-1}}-\dfrac{1}{\dbinom{a-1}{k-1}}=0$, we
conclude that $\dfrac{k-1}{k}\sum_{n=a}^{a-1}\dfrac{1}{\dbinom{n}{k}}%
=\dfrac{1}{\dbinom{a-1}{k-1}}-\dfrac{1}{\dbinom{a-1}{k-1}}$. In other words,
(\ref{pf.prop.binom.scaryfrac.goal}) holds for $b=a-1$. This completes the
induction base.

\textit{Induction step:} Let $\beta\in\left\{  a,a+1,a+2,\ldots\right\}  $.
Assume that (\ref{pf.prop.binom.scaryfrac.goal}) holds for $b=\beta-1$. We
must prove that (\ref{pf.prop.binom.scaryfrac.goal}) holds for $b=\beta$.

We have assumed that (\ref{pf.prop.binom.scaryfrac.goal}) holds for
$b=\beta-1$. In other words, we have%
\begin{equation}
\dfrac{k-1}{k}\sum_{n=a}^{\beta-1}\dfrac{1}{\dbinom{n}{k}}=\dfrac{1}%
{\dbinom{a-1}{k-1}}-\dfrac{1}{\dbinom{\beta-1}{k-1}}.
\label{pf.prop.binom.scaryfrac.IH}%
\end{equation}
(In particular, all fractions appearing in this equality are well-defined.)

We have $\beta\geq a$ (since $\beta\in\left\{  a,a+1,a+2,\ldots\right\}  $).

We have $k\in\left\{  1,2,3,\ldots\right\}  $ (since $k$ is a positive
integer). Hence, Proposition \ref{prop.binom.X-1} (applied to $\beta$ and $k$
instead of $m$ and $n$) yields $\dbinom{\beta}{k}=\dfrac{\beta}{k}%
\dbinom{\beta-1}{k-1}$. Multiplying this equality by $k$, we obtain%
\begin{equation}
k\dbinom{\beta}{k}=k\cdot\dfrac{\beta}{k}\dbinom{\beta-1}{k-1}=\beta
\dbinom{\beta-1}{k-1}. \label{pf.prop.binom.scaryfrac.4}%
\end{equation}
Now,%
\begin{align}
&  \dfrac{k-1}{k}\underbrace{\sum_{n=a}^{\beta}\dfrac{1}{\dbinom{n}{k}}%
}_{\substack{=\sum_{n=a}^{\beta-1}\dfrac{1}{\dbinom{n}{k}}+\dfrac{1}%
{\dbinom{\beta}{k}}\\\text{(since }\beta\geq a\text{)}}}\nonumber\\
&  =\dfrac{k-1}{k}\left(  \sum_{n=a}^{\beta-1}\dfrac{1}{\dbinom{n}{k}}%
+\dfrac{1}{\dbinom{\beta}{k}}\right) \nonumber\\
&  =\underbrace{\dfrac{k-1}{k}\sum_{n=a}^{\beta-1}\dfrac{1}{\dbinom{n}{k}}%
}_{\substack{=\dfrac{1}{\dbinom{a-1}{k-1}}-\dfrac{1}{\dbinom{\beta-1}{k-1}%
}\\\text{(by (\ref{pf.prop.binom.scaryfrac.IH}))}}}+\underbrace{\dfrac{k-1}%
{k}\cdot\dfrac{1}{\dbinom{\beta}{k}}}_{\substack{=\dfrac{k-1}{k\dbinom{\beta
}{k}}=\dfrac{k-1}{\beta\dbinom{\beta-1}{k-1}}\\\text{(by
(\ref{pf.prop.binom.scaryfrac.4}))}}}\nonumber\\
&  =\dfrac{1}{\dbinom{a-1}{k-1}}-\dfrac{1}{\dbinom{\beta-1}{k-1}}+\dfrac
{k-1}{\beta\dbinom{\beta-1}{k-1}}=\dfrac{1}{\dbinom{a-1}{k-1}}%
-\underbrace{\left(  \dfrac{1}{\dbinom{\beta-1}{k-1}}-\dfrac{k-1}{\beta
\dbinom{\beta-1}{k-1}}\right)  }_{\substack{=\dfrac{\beta-\left(  k-1\right)
}{\beta\dbinom{\beta-1}{k-1}}=\dfrac{\beta-k+1}{\beta\dbinom{\beta-1}{k-1}%
}\\\text{(since }\beta-\left(  k-1\right)  =\beta-k+1\text{)}}}\nonumber\\
&  =\dfrac{1}{\dbinom{a-1}{k-1}}-\dfrac{\beta-k+1}{\beta\dbinom{\beta-1}{k-1}%
}. \label{pf.prop.binom.scaryfrac.5}%
\end{align}

Also, Lemma \ref{lem.sol.binom.scaryfrac.absorp2} (applied to $m=\beta$ and
$n=k$) yields $\dbinom{\beta}{k}=\dfrac{\beta-k+1}{k}\dbinom{\beta}{k-1}$.
Multiplying this equality by $k$, we obtain%
\[
k\dbinom{\beta}{k}=k\cdot\dfrac{\beta-k+1}{k}\dbinom{\beta}{k-1}=\left(
\beta-k+1\right)  \dbinom{\beta}{k-1}.
\]
Comparing this with (\ref{pf.prop.binom.scaryfrac.4}), we obtain%
\[
\beta\dbinom{\beta-1}{k-1}=\left(  \beta-k+1\right)  \dbinom{\beta}{k-1}.
\]
Therefore,
\[
\dfrac{\beta-k+1}{\beta\dbinom{\beta-1}{k-1}}=\dfrac{\beta-k+1}{\left(
\beta-k+1\right)  \dbinom{\beta}{k-1}}=\dfrac{1}{\dbinom{\beta}{k-1}}.
\]
Hence, (\ref{pf.prop.binom.scaryfrac.5}) becomes%
\[
\dfrac{k-1}{k}\sum_{n=a}^{\beta}\dfrac{1}{\dbinom{n}{k}}=\dfrac{1}%
{\dbinom{a-1}{k-1}}-\underbrace{\dfrac{\beta-k+1}{\beta\dbinom{\beta-1}{k-1}}%
}_{=\dfrac{1}{\dbinom{\beta}{k-1}}}=\dfrac{1}{\dbinom{a-1}{k-1}}-\dfrac
{1}{\dbinom{\beta}{k-1}}.
\]
In other words, (\ref{pf.prop.binom.scaryfrac.goal}) holds for $b=\beta$. This
completes the induction step. Thus, (\ref{pf.prop.binom.scaryfrac.goal}) is
proven by induction. This completes the proof of Proposition
\ref{prop.binom.scaryfrac}.
\end{proof}

\begin{proof}
[Solution to Exercise \ref{exe.binom.scaryfrac}.]From $b\geq a\geq a-1$, we
obtain $b\in\left\{  a-1,a,a+1,\ldots\right\}  $ (since $b$ is an integer).
Thus, Proposition \ref{prop.binom.scaryfrac} yields $\dfrac{k-1}{k}\sum
_{n=a}^{b}\dfrac{1}{\dbinom{n}{k}}=\dfrac{1}{\dbinom{a-1}{k-1}}-\dfrac
{1}{\dbinom{b}{k-1}}$. This solves Exercise \ref{exe.binom.scaryfrac}.
\end{proof}

\subsection{Solution to Exercise \ref{exe.ps1.1.2}}

\begin{proof}
[Solution to Exercise \ref{exe.ps1.1.2}.]For every $N\in\mathbb{N}$, we let
$\left[  N\right]  $ denote the $N$-element set \newline$\left\{
1,2,\ldots,N\right\}  $.

For every $i\in\mathbb{N}$ and $j\in\mathbb{N}$, we define a \textit{filled
}$\left(  i,j\right)  $\textit{-set} to mean a subset $S$ of $\left[
i\right]  \times\left[  j\right]  $ satisfying the following two conditions:

\begin{enumerate}
\item For every $k\in\left[  i\right]  $, at least one element of $S$ has its
first coordinate\footnote{The \textit{coordinates} of a pair $\left(
u,v\right)  $ mean the entries $u$ and $v$. Thus, the first coordinate of
$\left(  u,v\right)  $ is $u$.} equal to $k$.

\item For every $\ell\in\left[  j\right]  $, at least one element of $S$ has
its second coordinate equal to $\ell$.
\end{enumerate}

We can visualize subsets $S$ of $\left[  i\right]  \times\left[  j\right]  $
as selections of boxes in a rectangular table that has $i$ rows and $j$
columns\footnote{Namely, for every $\left(  u,v\right)  \in S$, we select the
box in row $u$ and column $v$.}. For instance, for $i=3$ and $j=4$, we can
represent the subset $S=\left\{  \left(  1,1\right)  ,\left(  1,3\right)
,\left(  2,2\right)  ,\left(  3,1\right)  ,\left(  3,3\right)  ,\left(
3,4\right)  \right\}  $ of $\left[  i\right]  \times\left[  j\right]  $ as the
selection%
\begin{equation}%
\begin{tabular}
[c]{|c|c|c|c|}\hline
X &  & X & \\\hline
& X &  & \\\hline
X &  & X & X\\\hline
\end{tabular}
\ \ \label{sol.ps1.1.2.exa1}%
\end{equation}
(where the rows are labelled $1,2,3$ from top to bottom, the columns are
labelled $1,2,3,4$ from left to right, as in a matrix, and where the elements
of $S$ are marked with X'es). Condition 1 then says that every row contains at
least one selected box (i.e., at least one X); and Condition 2 says that every
column contains at least one selected box. Our example (\ref{sol.ps1.1.2.exa1}%
) satisfies these two conditions, but (for instance) the subset%
\[%
\begin{tabular}
[c]{|c|c|c|c|}\hline
X & X &  & X\\\hline
& X & \phantom{X} & \\\hline
X &  &  & X\\\hline
\end{tabular}
\ \
\]
does not (it fails Condition 2).

For every $i\in\mathbb{N}$ and $j\in\mathbb{N}$, we define a nonnegative
integer $c_{i,j}$ as the number of all filled $\left(  i,j\right)  $-sets
which have $n$ elements\footnote{When we say \textquotedblleft have $n$
elements\textquotedblright, we mean \textquotedblleft have exactly $n$
elements\textquotedblright, not \textquotedblleft have at least $n$
elements\textquotedblright.}. We shall now show that (\ref{eq.exe.1.2.claim})
is satisfied.

Indeed, let us first prove that any $x\in\mathbb{N}$ and $y\in\mathbb{N}$
satisfy%
\begin{equation}
\dbinom{xy}{n}=\sum_{i=0}^{n}\sum_{j=0}^{n}c_{i,j}\dbinom{x}{i}\dbinom{y}{j}.
\label{sol.ps1.1.2.xyclaim}%
\end{equation}

Keep in mind that (\ref{sol.ps1.1.2.xyclaim}) and (\ref{eq.exe.1.2.claim}) are
different claims: The $x$ and $y$ in (\ref{sol.ps1.1.2.xyclaim}) are
nonnegative integers, while the $X$ and $Y$ in (\ref{eq.exe.1.2.claim}) are indeterminates!

[\textit{Proof of (\ref{sol.ps1.1.2.xyclaim}):} Let $x\in\mathbb{N}$ and
$y\in\mathbb{N}$. Recall that $\dbinom{xy}{n}$ is the number of all
$n$-element subsets of a given $xy$-element set. Since $\left[  x\right]
\times\left[  y\right]  $ is an $xy$-element set, we thus conclude that
$\dbinom{xy}{n}$ is the number of all $n$-element subsets of $\left[
x\right]  \times\left[  y\right]  $.

Now, let us find a different way to count all $n$-element subsets of $\left[
x\right]  \times\left[  y\right]  $. As above, we can visualize such subsets
as selections of boxes in a rectangular table that has $x$ rows and $y$
columns; we again mark the selected boxes by X'es. We want to count all ways
to select $n$ boxes in this table, i.e., to place $n$ X'es in the table. We
can place $n$ X'es in the table by means of the following process:

\begin{enumerate}
\item We choose how many rows of the table will have at least one X. This can
be a number from $0$ to $n$ (inclusive)\footnote{I am not saying that any
number from $0$ to $n$ (inclusive) is possible; I am just saying that this
will always be a number from $0$ to $n$ (inclusive). Here is why:
\par
Clearly, the number of rows of the table that will have at least one X is
$\geq0$. But it is also $\leq n$, because we want to place only $n$ X'es in
the table, and these $n$ X'es will clearly occupy at most $n$ rows.}; we
denote it by $i$.

\item We choose how many columns of the table will have at least one X. This
can be a number from $0$ to $n$ (inclusive)\footnote{This follows by a similar
argument as the analogous statement in Step 1.}; we denote it by $j$.

\item We choose the $i$ rows of the table that will have at least one X. This
can be done in $\dbinom{x}{i}$ ways (since there are $x$ rows to choose from).

\item We choose the $j$ columns of the table that will have at least one X.
This can be done in $\dbinom{y}{j}$ ways (since there are $y$ columns to
choose from).

\item It remains to place $n$ X'es in the table in such a way that the rows
that contain at least one X are precisely the $i$ chosen rows, and the columns
that contain at least one X are precisely the $j$ chosen columns. To do so, we
can temporarily remove all the remaining $x-i$ rows and $y-j$ columns. We are
then left with a rectangular table that has $i$ rows and $j$ columns, and now
we need to place $n$ X'es in it in such a way that every row contains at least
one X and every column contains at least one X. As we know, the number of ways
to do this is $c_{i,j}$ (because this is how $c_{i,j}$ was defined).
\end{enumerate}

This process makes it clear that the total number of ways to place $n$ X'es in
the (original) table is $\sum_{i=0}^{n}\sum_{j=0}^{n}\dbinom{x}{i}\dbinom
{y}{j}c_{i,j}=\sum_{i=0}^{n}\sum_{j=0}^{n}c_{i,j}\dbinom{x}{i}\dbinom{y}{j}$.
In other words, the number of all $n$-element subsets of $\left[  x\right]
\times\left[  y\right]  $ is $\sum_{i=0}^{n}\sum_{j=0}^{n}c_{i,j}\dbinom{x}%
{i}\dbinom{y}{j}$.

So we know that this number equals both $\dbinom{xy}{n}$ and $\sum_{i=0}%
^{n}\sum_{j=0}^{n}c_{i,j}\dbinom{x}{i}\dbinom{y}{j}$ at the same time.
Comparing these values, we obtain (\ref{sol.ps1.1.2.xyclaim}).]

Now that (\ref{sol.ps1.1.2.xyclaim}) is proven, we can finally solve the
exercise. We define two polynomials $P$ and $Q$ in the indeterminates $X$ and
$Y$ with rational coefficients by setting%
\begin{align*}
P  &  =\dbinom{XY}{n};\\
Q  &  =\sum_{i=0}^{n}\sum_{j=0}^{n}c_{i,j}\dbinom{X}{i}\dbinom{Y}{j}%
\end{align*}
\footnote{These are both polynomials since $\dbinom{XY}{n}$, $\dbinom{X}{i}$
and $\dbinom{Y}{j}$ are polynomials in $X$ and $Y$.}. The equality
(\ref{sol.ps1.1.2.xyclaim}) (which we have proven) states that $P\left(
x,y\right)  =Q\left(  x,y\right)  $ for all $x\in\mathbb{N}$ and
$y\in\mathbb{N}$. Thus, Lemma \ref{lem.polyeq} \textbf{(d)} yields that $P=Q$.
Recalling how $P$ and $Q$ are defined, we see that this is precisely the
equality (\ref{eq.exe.1.2.claim}).

Hence, (\ref{eq.exe.1.2.claim}) is proven, and Exercise \ref{exe.ps1.1.2} solved.
\end{proof}

\begin{remark}
I learnt the above solution to Exercise \ref{exe.ps1.1.2}
\href{http://www.artofproblemsolving.com/community/c6h299793p1623722}{from
Gjergji Zaimi on AoPS}. The numbers $c_{i,j}$ constructed in the solution do
not appear to be easily computable by a simple closed formula. Nevertheless,
they have some nice properties that can be easily obtained from their
combinatorial definition:

\begin{itemize}
\item We have $c_{i,j}=c_{j,i}$ for all $i\in\mathbb{N}$ and $j\in\mathbb{N}$.

\item We have $c_{0,0}=%
\begin{cases}
1, & \text{if }n=0;\\
0, & \text{if }n>0
\end{cases}
$, but every positive integer $j$ satisfies $c_{0,j}=0$.

\item We have $c_{1,j}=%
\begin{cases}
1, & \text{if }n=j;\\
0, & \text{if }n\neq j
\end{cases}
$ for all positive integers $j$.

\item We have $c_{2,2}=%
\begin{cases}
0, & \text{if }n\leq1;\\
2, & \text{if }n=2;\\
4, & \text{if }n=3;\\
1, & \text{if }n=4;\\
0, & \text{if }n>4
\end{cases}
$.

\item We have $c_{i,j}=0$ if $n>ij$.

\item For every $j\in\mathbb{N}$, the number $c_{n,j}$ is the number of all
surjective maps $\left\{  1,2,\ldots,n\right\}  \rightarrow\left\{
1,2,\ldots,j\right\}  $.
\end{itemize}
\end{remark}

\subsection{Solution to Exercise \ref{exe.ps1.1.3}}

\begin{proof}
[Solution to Exercise \ref{exe.ps1.1.3}.]Here is one possible solution:

Exercise \ref{exe.ps1.1.2} shows that, for every $n\in\mathbb{N}$, there exist
\textbf{nonnegative} integers $c_{i,j}$ for all $0\leq i\leq n$ and $0\leq
j\leq n$ such that (\ref{eq.exe.1.2.claim}) holds. We denote these integers
$c_{i,j}$ by $c_{i,j,n}$ (in order to make their dependence on $n$ explicit).
Thus, for every $n\in\mathbb{N}$, the nonnegative integers $c_{i,j,n}$ defined
for all $0\leq i\leq n$ and $0\leq j\leq n$ satisfy%
\begin{equation}
\dbinom{XY}{n}=\sum_{i=0}^{n}\sum_{j=0}^{n}c_{i,j,n}\dbinom{X}{i}\dbinom{Y}%
{j}. \label{pf.exe.1.3.1}%
\end{equation}
Substituting $a$ and $X$ for $X$ and $Y$ in this equality, we obtain%
\begin{equation}
\dbinom{aX}{n}=\sum_{i=0}^{n}\sum_{j=0}^{n}c_{i,j,n}\dbinom{a}{i}\dbinom{X}%
{j}. \label{pf.exe.1.3.2}%
\end{equation}

Now, Theorem \ref{thm.vandermonde.XY} (applied to $n=c$) yields%
\[
\dbinom{X+Y}{c}=\sum_{k=0}^{c}\dbinom{X}{k}\dbinom{Y}{c-k}.
\]
Substituting $aX$ and $b$ for $X$ and $Y$ in this equality, we obtain%
\begin{align}
\dbinom{aX+b}{c}  &  =\sum_{k=0}^{c}\underbrace{\dbinom{aX}{k}}%
_{\substack{=\sum_{i=0}^{k}\sum_{j=0}^{k}c_{i,j,k}\dbinom{a}{i}\dbinom{X}%
{j}\\\text{(by (\ref{pf.exe.1.3.2}), applied to }n=k\text{)}}}\dbinom{b}%
{c-k}\nonumber\\
&  =\sum_{k=0}^{c}\left(  \sum_{i=0}^{k}\sum_{j=0}^{k}c_{i,j,k}\dbinom{a}%
{i}\dbinom{X}{j}\right)  \dbinom{b}{c-k}\nonumber\\
&  =\underbrace{\sum_{k=0}^{c}\sum_{i=0}^{k}\sum_{j=0}^{k}}_{=\sum_{j=0}%
^{c}\sum_{k=j}^{c}\sum_{i=0}^{k}}c_{i,j,k}\dbinom{a}{i}\underbrace{\dbinom
{X}{j}\dbinom{b}{c-k}}_{=\dbinom{b}{c-k}\dbinom{X}{j}}\nonumber\\
&  =\sum_{j=0}^{c}\sum_{k=j}^{c}\sum_{i=0}^{k}c_{i,j,k}\dbinom{a}{i}\dbinom
{b}{c-k}\dbinom{X}{j}. \label{pf.exe.1.3.5}%
\end{align}

Now, for every $j\in\left\{  0,1,\ldots,c\right\}  $, define an integer
$d_{j}$ by $d_{j}=\sum_{k=j}^{c}\sum_{i=0}^{k}c_{i,j,k}\dbinom{a}{i}\dbinom
{b}{c-k}$. This $d_{j}$ is clearly a nonnegative integer (since the
$c_{i,j,n}$ are nonnegative, and so are the binomial coefficients $\dbinom
{a}{i}$ and $\dbinom{b}{c-k}$ (due to $a$ and $b$ being nonnegative)). Then,
(\ref{pf.exe.1.3.5}) becomes%
\[
\dbinom{aX+b}{c}=\sum_{j=0}^{c}\underbrace{\sum_{k=j}^{c}\sum_{i=0}%
^{k}c_{i,j,k}\dbinom{a}{i}\dbinom{b}{c-k}}_{=d_{j}}\dbinom{X}{j}=\sum
_{j=0}^{c}d_{j}\dbinom{X}{j}=\sum_{i=0}^{c}d_{i}\dbinom{X}{i}.
\]
Exercise \ref{exe.ps1.1.3} is thus solved.
\end{proof}

\subsection{Solution to Exercise \ref{exe.binom.rewr-prod}}

Before we come to the solution of Exercise \ref{exe.binom.rewr-prod}, let us
show a simple identity:

\begin{proposition}
\label{prop.pascal-rowsum+0}Let $m\in\mathbb{N}$ and $p\in\mathbb{N}$ be such
that $p\geq m$. Then,%
\[
\sum_{k=0}^{p}\dbinom{m}{k}=2^{m}.
\]

\end{proposition}

\begin{proof}
[Proof of Proposition \ref{prop.pascal-rowsum+0}.]For each $k\in\left\{
m+1,m+2,\ldots,p\right\}  $, we have%
\begin{equation}
\dbinom{m}{k}=0 \label{pf.prop.pascal-rowsum+0.1}%
\end{equation}
\footnote{\textit{Proof of (\ref{pf.prop.pascal-rowsum+0.1}):} Let
$k\in\left\{  m+1,m+2,\ldots,p\right\}  $. Thus, $k\geq m+1>m$ and therefore
$m<k$. Also, $k>m\geq0$ (since $m\in\mathbb{N}$). Hence, $k\in\mathbb{N}$.
Therefore, Proposition \ref{prop.binom.0} (applied to $n=k$) yields
$\dbinom{m}{k}=0$. This proves (\ref{pf.prop.pascal-rowsum+0.1}).}.

We have $p\geq m\geq0$ (since $m\in\mathbb{N}$). Thus, the sum $\sum_{k=0}%
^{p}\dbinom{m}{k}$ can be split as follows:%
\[
\sum_{k=0}^{p}\dbinom{m}{k}=\sum_{k=0}^{m}\dbinom{m}{k}+\sum_{k=m+1}%
^{p}\underbrace{\dbinom{m}{k}}_{\substack{=0\\\text{(by
(\ref{pf.prop.pascal-rowsum+0.1}))}}}=\sum_{k=0}^{m}\dbinom{m}{k}%
+\underbrace{\sum_{k=m+1}^{p}0}_{=0}=\sum_{k=0}^{m}\dbinom{m}{k}=2^{m}%
\]
(by Proposition \ref{prop.binom.bin-id} \textbf{(b)} (applied to $n=m$)). This
proves Proposition \ref{prop.pascal-rowsum+0}.
\end{proof}

\begin{proof}
[Solution to Exercise \ref{exe.binom.rewr-prod}.]\textbf{(a)} For each
$k\in\mathbb{N}$, we have%
\begin{align}
&  \underbrace{\dbinom{n}{k}}_{\substack{=\dfrac{n\left(  n-1\right)
\cdots\left(  n-k+1\right)  }{k!}\\\text{(by (\ref{eq.binom.mn}) (applied to
}n\text{ and }k\\\text{instead of }m\text{ and }n\text{))}}%
}\underbrace{\dbinom{n-k}{b}}_{\substack{=\dfrac{\left(  n-k\right)  \left(
\left(  n-k\right)  -1\right)  \cdots\left(  \left(  n-k\right)  -b+1\right)
}{b!}\\\text{(by (\ref{eq.binom.mn}) (applied to }n-k\text{ and }%
b\\\text{instead of }m\text{ and }n\text{))}}}\nonumber\\
&  =\dfrac{n\left(  n-1\right)  \cdots\left(  n-k+1\right)  }{k!}\cdot
\dfrac{\left(  n-k\right)  \left(  \left(  n-k\right)  -1\right)
\cdots\left(  \left(  n-k\right)  -b+1\right)  }{b!}\nonumber\\
&  =\dfrac{1}{k!b!}\cdot\left(  n\left(  n-1\right)  \cdots\left(
n-k+1\right)  \right)  \cdot\left(  \underbrace{\left(  n-k\right)  \left(
\left(  n-k\right)  -1\right)  \cdots\left(  \left(  n-k\right)  -b+1\right)
}_{=\left(  n-k\right)  \left(  n-k-1\right)  \cdots\left(  n-k-b+1\right)
}\right) \nonumber\\
&  =\dfrac{1}{k!b!}\cdot\underbrace{\left(  n\left(  n-1\right)  \cdots\left(
n-k+1\right)  \right)  \cdot\left(  \left(  n-k\right)  \left(  n-k-1\right)
\cdots\left(  n-k-b+1\right)  \right)  }_{=n\left(  n-1\right)  \cdots\left(
n-k-b+1\right)  }\nonumber\\
&  =\dfrac{1}{k!b!}\cdot\left(  n\left(  n-1\right)  \cdots\left(
n-k-b+1\right)  \right)  . \label{sol.binom.rewr-prod.a.0}%
\end{align}

Let $j\geq a$ be an integer. Thus, $j \geq a \geq0$ (since $a \in\mathbb{N}$),
so that $j \in\mathbb{N}$. Hence, Proposition \ref{prop.binom.formula}
(applied to $j$ and $a$ instead of $m$ and $n$) shows that%
\[
\dbinom{j}{a}=\dfrac{j!}{a!\left(  j-a\right)  !}.
\]
But (\ref{sol.binom.rewr-prod.a.0}) (applied to $k=j$) shows that%
\[
\dbinom{n}{j}\dbinom{n-j}{b}=\dfrac{1}{j!b!}\cdot\left(  n\left(  n-1\right)
\cdots\left(  n-j-b+1\right)  \right)  .
\]
Multiplying these two equalities, we obtain%
\begin{align}
\dbinom{j}{a}\dbinom{n}{j}\dbinom{n-j}{b}  &  =\dfrac{j!}{a!\left(
j-a\right)  !}\cdot\dfrac{1}{j!b!}\cdot\left(  n\left(  n-1\right)
\cdots\left(  n-j-b+1\right)  \right) \nonumber\\
&  =\dfrac{1}{a!\left(  j-a\right)  !b!}\cdot\left(  n\left(  n-1\right)
\cdots\left(  n-j-b+1\right)  \right)  . \label{sol.binom.rewr-prod.a.1}%
\end{align}

On the other hand, set $m=n-a-b$. Then,%
\[
\underbrace{m}_{=n-a-b}-\left(  j-a\right)  =\left(  n-a-b\right)  -\left(
j-a\right)  =n-j-b
\]
and $n\geq m$ (since $m=n-\underbrace{a}_{\geq0}-\underbrace{b}_{\geq0}\leq
n$) and $m\geq m-\left(  j-a\right)  $ (since $m-\left(  \underbrace{j}_{\geq
a}-a\right)  \leq m-\left(  a-a\right)  =m$).

The equality (\ref{sol.binom.rewr-prod.a.0}) (applied to $k=a$) yields%
\[
\dbinom{n}{a}\dbinom{n-a}{b}=\dfrac{1}{a!b!}\cdot\left(  n\left(  n-1\right)
\cdots\left(  n-a-b+1\right)  \right)  =\dfrac{1}{a!b!}\cdot\left(  n\left(
n-1\right)  \cdots\left(  m+1\right)  \right)
\]
(since $n-a-b=m$). Furthermore, $j-a\in\mathbb{N}$ (since $j\geq a$), and thus
the binomial coefficient $\dbinom{m}{j-a}$ is well-defined. Furthermore, from
$n-a-b=m$, we obtain%
\[
\dbinom{n-a-b}{j-a}=\dbinom{m}{j-a}=\dfrac{m\left(  m-1\right)  \cdots\left(
m-\left(  j-a\right)  +1\right)  }{\left(  j-a\right)  !}%
\]
(by (\ref{eq.binom.mn}) (applied to $j-a$ instead of $n$)). Now,
\begin{align*}
&  \underbrace{\dbinom{n}{a}\dbinom{n-a}{b}}_{=\dfrac{1}{a!b!}\cdot\left(
n\left(  n-1\right)  \cdots\left(  m+1\right)  \right)  }%
\ \ \ \underbrace{\dbinom{n-a-b}{j-a}}_{=\dfrac{m\left(  m-1\right)
\cdots\left(  m-\left(  j-a\right)  +1\right)  }{\left(  j-a\right)  !}}\\
&  =\dfrac{1}{a!b!}\cdot\left(  n\left(  n-1\right)  \cdots\left(  m+1\right)
\right)  \cdot\dfrac{m\left(  m-1\right)  \cdots\left(  m-\left(  j-a\right)
+1\right)  }{\left(  j-a\right)  !}\\
&  =\dfrac{1}{a!\left(  j-a\right)  !b!}\cdot\underbrace{\left(  n\left(
n-1\right)  \cdots\left(  m+1\right)  \right)  \cdot\left(  m\left(
m-1\right)  \cdots\left(  m-\left(  j-a\right)  +1\right)  \right)
}_{\substack{=n\left(  n-1\right)  \cdots\left(  m-\left(  j-a\right)
+1\right)  \\\text{(since }n\geq m\geq m-\left(  j-a\right)  \text{)}}}\\
&  =\dfrac{1}{a!\left(  j-a\right)  !b!}\cdot\left(  n\left(  n-1\right)
\cdots\left(  m-\left(  j-a\right)  +1\right)  \right) \\
&  =\dfrac{1}{a!\left(  j-a\right)  !b!}\cdot\left(  n\left(  n-1\right)
\cdots\left(  n-j-b+1\right)  \right)  \ \ \ \ \ \ \ \ \ \ \left(  \text{since
}m-\left(  j-a\right)  =n-j-b\right)  .
\end{align*}
Comparing this with (\ref{sol.binom.rewr-prod.a.1}), we obtain
\[
\dbinom{n}{a}\dbinom{n-a}{b}\dbinom{n-a-b}{j-a}=\dbinom{j}{a}\dbinom{n}%
{j}\dbinom{n-j}{b}=\dbinom{n}{j}\dbinom{j}{a}\dbinom{n-j}{b}.
\]
This solves Exercise \ref{exe.binom.rewr-prod} \textbf{(a)}.

\textbf{(b)} Let $n\geq a$ be an integer. Thus, $n-a\in\mathbb{N}$. Now, we
claim that%
\begin{equation}
\sum_{j=a}^{n}\dbinom{n}{j}\dbinom{j}{a}\dbinom{n-j}{b}=\dbinom{n}{a}%
\dbinom{n-a}{b}2^{n-a-b}. \label{sol.binom.rewr-prod.b.claim}%
\end{equation}

[\textit{Proof of (\ref{sol.binom.rewr-prod.b.claim}):} First of all,
(\ref{sol.binom.rewr-prod.b.claim}) holds if $n-a-b<0$%
\ \ \ \ \footnote{\textit{Proof.} Assume that $n-a-b<0$. Thus, $n-a<b$. Hence,
Proposition \ref{prop.binom.0} (applied to $n-a$ and $b$ instead of $m$ and
$n$) shows that $\dbinom{n-a}{b}=0$ (since $n-a\in\mathbb{N}$). But%
\begin{align*}
\sum_{j=a}^{n}\underbrace{\dbinom{n}{j}\dbinom{j}{a}\dbinom{n-j}{b}%
}_{\substack{=\dbinom{n}{a}\dbinom{n-a}{b}\dbinom{n-a-b}{j-a}\\\text{(by
Exercise \ref{exe.binom.rewr-prod} \textbf{(a)})}}}  &  =\sum_{j=a}^{n}%
\dbinom{n}{a}\underbrace{\dbinom{n-a}{b}}_{=0}\dbinom{n-a-b}{j-a}\\
&  =\sum_{j=a}^{n}\dbinom{n}{a}0\dbinom{n-a-b}{j-a}=0.
\end{align*}
Comparing this with%
\[
\dbinom{n}{a}\underbrace{\dbinom{n-a}{b}}_{=0}2^{n-a-b}=0,
\]
we obtain $\sum_{j=a}^{n}\dbinom{n}{j}\dbinom{j}{a}\dbinom{n-j}{b}=\dbinom
{n}{a}\dbinom{n-a}{b}2^{n-a-b}$. Hence, (\ref{sol.binom.rewr-prod.b.claim}) is
proven under the assumption that $n-a-b<0$.}. Hence, for the rest of this
proof, we can WLOG assume that we don't have $n-a-b<0$. Assume this.

We have $n-a-b\geq0$ (since we don't have $n-a-b<0$), so that $n-a-b\in
\mathbb{N}$. Thus, we have $n-a\in\mathbb{N}$ and $n-a-b\in\mathbb{N}$ and
$n-a\geq n-a-b$ (since $\left(  n-a\right)  -\left(  n-a-b\right)  =b\geq0$).
Hence, Proposition \ref{prop.pascal-rowsum+0} (applied to $m=n-a-b$ and
$p=n-a$) yields%
\[
\sum_{k=0}^{n-a}\dbinom{n-a-b}{k}=2^{n-a-b}.
\]

Now,
\begin{align*}
&  \sum_{j=a}^{n}\underbrace{\dbinom{n}{j}\dbinom{j}{a}\dbinom{n-j}{b}%
}_{\substack{=\dbinom{n}{a}\dbinom{n-a}{b}\dbinom{n-a-b}{j-a}\\\text{(by
Exercise \ref{exe.binom.rewr-prod} \textbf{(a)})}}}\\
&  =\sum_{j=a}^{n}\dbinom{n}{a}\dbinom{n-a}{b}\dbinom{n-a-b}{j-a}=\dbinom
{n}{a}\dbinom{n-a}{b}\sum_{j=a}^{n}\dbinom{n-a-b}{j-a}\\
&  =\dbinom{n}{a}\dbinom{n-a}{b}\underbrace{\sum_{k=0}^{n-a}\dbinom{n-a-b}{k}%
}_{=2^{n-a-b}}\\
&  \ \ \ \ \ \ \ \ \ \ \left(  \text{here, we have substituted }k\text{ for
}j-a\text{ in the sum}\right) \\
&  =\dbinom{n}{a}\dbinom{n-a}{b}2^{n-a-b}.
\end{align*}
This proves (\ref{sol.binom.rewr-prod.b.claim}).]

Clearly, (\ref{sol.binom.rewr-prod.b.claim}) answers Exercise
\ref{exe.binom.rewr-prod} \textbf{(b)}.
\end{proof}

\subsection{Solution to Exercise \ref{exe.pie.binid}}

There are many ways to solve Exercise \ref{exe.pie.binid}. The following one
might be the shortest:

\begin{proof}
[Proof of Lemma \ref{lem.pie.binid}.]We proceed by induction over $m$:

\begin{vershort}
\textit{Induction base:} We have
\begin{align*}
&  \sum_{r=0}^{k}\left(  -1\right)  ^{r}\dbinom{n}{r}\dbinom{r}{k}\\
&  =\sum_{r=0}^{k-1}\left(  -1\right)  ^{r}\dbinom{n}{r}\underbrace{\dbinom
{r}{k}}_{\substack{=0\\\text{(by Proposition \ref{prop.binom.0} (applied to
}r\text{ and }k\\\text{instead of }m\text{ and }n\text{) (since }r\leq
k-1<k\text{))}}}+\left(  -1\right)  ^{k}\dbinom{n}{k}\underbrace{\dbinom{k}%
{k}}_{\substack{=1\\\text{(by Proposition \ref{prop.binom.mm}}\\\text{(applied
to }k\text{ instead of }m\text{))}}}\\
&  \ \ \ \ \ \ \ \ \ \ \left(  \text{here, we have split off the addend for
}r=k\text{ from the sum}\right) \\
&  =\underbrace{\sum_{r=0}^{k-1}\left(  -1\right)  ^{r}\dbinom{n}{r}0}%
_{=0}+\left(  -1\right)  ^{k}\dbinom{n}{k}=\left(  -1\right)  ^{k}\dbinom
{n}{k}.
\end{align*}
Comparing this with%
\[
\left(  -1\right)  ^{k}\dbinom{n}{k}\underbrace{\dbinom{n-k-1}{k-k}%
}_{\substack{=\dbinom{n-k-1}{0}=1\\\text{(by Proposition \ref{prop.binom.00}
\textbf{(a)}}\\\text{(applied to }n-k-1\text{ instead of }m\text{))}}}=\left(
-1\right)  ^{k}\dbinom{n}{k},
\]
we obtain $\sum_{r=0}^{k}\left(  -1\right)  ^{r}\dbinom{n}{r}\dbinom{r}%
{k}=\left(  -1\right)  ^{k}\dbinom{n}{k}\dbinom{n-k-1}{k-k}$. In other words,
Lemma \ref{lem.pie.binid} holds for $m=k$. This completes the induction base.
\end{vershort}

\begin{verlong}
\textit{Induction base:} Assume that $m=k$. Notice that%
\begin{equation}
\dbinom{r}{k}=0 \label{sol.pie.binid.base.1}%
\end{equation}
for each $r\in\left\{  0,1,\ldots,k-1\right\}  $%
\ \ \ \ \footnote{\textit{Proof of (\ref{sol.pie.binid.base.1}):} Let
$r\in\left\{  0,1,\ldots,k-1\right\}  $. Thus, $0\leq r\leq k-1$. Also,
$r\in\left\{  0,1,\ldots,k-1\right\}  \subseteq\mathbb{N}$. We have $r\leq
k-1<k$. Hence, Proposition \ref{prop.binom.0} (applied to $r$ and $k$ instead
of $m$ and $n$) yields $\dbinom{r}{k}=0$. This proves
(\ref{sol.pie.binid.base.1}).}.

But from $m=k$, we obtain%
\begin{align}
\sum_{r=0}^{m}\left(  -1\right)  ^{r}\dbinom{n}{r}\dbinom{r}{k}  &
=\sum_{r=0}^{k}\left(  -1\right)  ^{r}\dbinom{n}{r}\dbinom{r}{k}\nonumber\\
&  =\sum_{r=0}^{k-1}\left(  -1\right)  ^{r}\dbinom{n}{r}\underbrace{\dbinom
{r}{k}}_{\substack{=0\\\text{(by (\ref{sol.pie.binid.base.1}))}}}+\left(
-1\right)  ^{k}\dbinom{n}{k}\underbrace{\dbinom{k}{k}}%
_{\substack{=1\\\text{(by Proposition \ref{prop.binom.mm}}\\\text{(applied to
}k\text{ instead of }m\text{))}}}\nonumber\\
&  \ \ \ \ \ \ \ \ \ \ \left(  \text{here, we have split off the addend for
}r=k\text{ from the sum}\right) \nonumber\\
&  =\underbrace{\sum_{r=0}^{k-1}\left(  -1\right)  ^{r}\dbinom{n}{r}0}%
_{=0}+\left(  -1\right)  ^{k}\dbinom{n}{k}=\left(  -1\right)  ^{k}\dbinom
{n}{k}. \label{sol.pie.binid.base.L}%
\end{align}

On the other hand, $m-k=0$ (since $m=k$), so that%
\[
\dbinom{n-k-1}{m-k}=\dbinom{n-k-1}{0}=1
\]
(by Proposition \ref{prop.binom.00} \textbf{(a)} (applied to $n-k-1$ instead
of $m$)). Hence,%
\[
\underbrace{\left(  -1\right)  ^{m}}_{\substack{=\left(  -1\right)
^{k}\\\text{(since }m=k\text{)}}}\dbinom{n}{k}\underbrace{\dbinom{n-k-1}{m-k}%
}_{=1}=\left(  -1\right)  ^{k}\dbinom{n}{k}.
\]
Comparing this with (\ref{sol.pie.binid.base.L}), we find
\[
\sum_{r=0}^{m}\left(  -1\right)  ^{r}\dbinom{n}{r}\dbinom{r}{k}=\left(
-1\right)  ^{m}\dbinom{n}{k}\dbinom{n-k-1}{m-k}.
\]
Hence, Lemma \ref{lem.pie.binid} is proven for $m=k$. This completes the
induction base.
\end{verlong}

\textit{Induction step:} Let $M$ be an element of $\left\{  k,k+1,k+2,\ldots
\right\}  $ such that $M>k$. Assume that Lemma \ref{lem.pie.binid} holds for
$m=M-1$. We must now prove that Lemma \ref{lem.pie.binid} holds for $m=M$.

We have assumed that Lemma \ref{lem.pie.binid} holds for $m=M-1$. In other
words, we have%
\[
\sum_{r=0}^{M-1}\left(  -1\right)  ^{r}\dbinom{n}{r}\dbinom{r}{k}=\left(
-1\right)  ^{M-1}\dbinom{n}{k}\dbinom{n-k-1}{\left(  M-1\right)  -k}.
\]
Hence,%
\begin{align}
\sum_{r=0}^{M-1}\left(  -1\right)  ^{r}\dbinom{n}{r}\dbinom{r}{k}  &
=\underbrace{\left(  -1\right)  ^{M-1}}_{=-\left(  -1\right)  ^{M}}\dbinom
{n}{k}\underbrace{\dbinom{n-k-1}{\left(  M-1\right)  -k}}_{\substack{=\dbinom
{n-k-1}{M-k-1}\\\text{(since }\left(  M-1\right)  -k=M-k-1\text{)}%
}}\nonumber\\
&  =-\left(  -1\right)  ^{M}\dbinom{n}{k}\dbinom{n-k-1}{M-k-1}.
\label{sol.pie.binid.step.hyp2}%
\end{align}

On the other hand, $M\geq k$ (since $M>k$). Hence, Proposition
\ref{prop.binom.trinom-rev} (applied to $n$, $M$ and $k$ instead of $m$, $i$
and $a$) yields%
\begin{equation}
\dbinom{n}{M}\dbinom{M}{k}=\dbinom{n}{k}\dbinom{n-k}{M-k}.
\label{sol.pie.binid.step.3}%
\end{equation}

Furthermore, $M-k>0$ (since $M>k$), so that $M-k\in\left\{  1,2,3,\ldots
\right\}  $ (since $M-k$ is an integer). Hence, Proposition
\ref{prop.binom.rec} (applied to $n-k$ and $M-k$ instead of $m$ and $n$)
yields%
\[
\dbinom{n-k}{M-k}=\dbinom{n-k-1}{M-k-1}+\dbinom{n-k-1}{M-k}.
\]
Hence, (\ref{sol.pie.binid.step.3}) becomes%
\begin{align}
\dbinom{n}{M}\dbinom{M}{k}  &  =\dbinom{n}{k}\underbrace{\dbinom{n-k}{M-k}%
}_{=\dbinom{n-k-1}{M-k-1}+\dbinom{n-k-1}{M-k}}\nonumber\\
&  =\dbinom{n}{k}\left(  \dbinom{n-k-1}{M-k-1}+\dbinom{n-k-1}{M-k}\right)  .
\label{sol.pie.binid.step.5}%
\end{align}

Now,
\begin{align*}
&  \sum_{r=0}^{M}\left(  -1\right)  ^{r}\dbinom{n}{r}\dbinom{r}{k}\\
&  =\underbrace{\sum_{r=0}^{M-1}\left(  -1\right)  ^{r}\dbinom{n}{r}\dbinom
{r}{k}}_{\substack{=-\left(  -1\right)  ^{M}\dbinom{n}{k}\dbinom{n-k-1}%
{M-k-1}\\\text{(by (\ref{sol.pie.binid.step.hyp2}))}}}+\left(  -1\right)
^{M}\underbrace{\dbinom{n}{M}\dbinom{M}{k}}_{\substack{=\dbinom{n}{k}\left(
\dbinom{n-k-1}{M-k-1}+\dbinom{n-k-1}{M-k}\right)  \\\text{(by
(\ref{sol.pie.binid.step.5}))}}}\\
&  \ \ \ \ \ \ \ \ \ \ \left(  \text{here, we have split off the addend for
}r=M\text{ from the sum}\right) \\
&  =-\left(  -1\right)  ^{M}\dbinom{n}{k}\dbinom{n-k-1}{M-k-1}%
+\underbrace{\left(  -1\right)  ^{M}\dbinom{n}{k}\left(  \dbinom{n-k-1}%
{M-k-1}+\dbinom{n-k-1}{M-k}\right)  }_{=\left(  -1\right)  ^{M}\dbinom{n}%
{k}\dbinom{n-k-1}{M-k-1}+\left(  -1\right)  ^{M}\dbinom{n}{k}\dbinom
{n-k-1}{M-k}}\\
&  =-\left(  -1\right)  ^{M}\dbinom{n}{k}\dbinom{n-k-1}{M-k-1}+\left(
-1\right)  ^{M}\dbinom{n}{k}\dbinom{n-k-1}{M-k-1}+\left(  -1\right)
^{M}\dbinom{n}{k}\dbinom{n-k-1}{M-k}\\
&  =\left(  -1\right)  ^{M}\dbinom{n}{k}\dbinom{n-k-1}{M-k}.
\end{align*}
In other words, Lemma \ref{lem.pie.binid} holds for $m=M$. This completes the
induction step. Hence, Lemma \ref{lem.pie.binid} is proven.
\end{proof}

We notice that Lemma \ref{lem.pie.binid} holds even if we replace the
assumption \textquotedblleft$n\in\mathbb{N}$\textquotedblright\ by
\textquotedblleft$n\in\mathbb{Q}$\textquotedblright. In fact, our above proof
of Lemma \ref{lem.pie.binid} applies verbatim in this more general setting.

\begin{proof}
[Solution to Exercise \ref{exe.pie.binid}.]We have proven Lemma
\ref{lem.pie.binid}; hence, Exercise \ref{exe.pie.binid} is solved.
\end{proof}

\subsection{Solution to Exercise \ref{exe.iverson-prop}}

\begin{proof}
[Solution to Exercise \ref{exe.iverson-prop}.]\textbf{(a)} Let $\mathcal{A}$
and $\mathcal{B}$ be two equivalent statements.

We have $\left[  \mathcal{A}\right]  =%
\begin{cases}
1, & \text{if }\mathcal{A}\text{ is true};\\
0, & \text{if }\mathcal{A}\text{ is false}%
\end{cases}
$ (by the definition of $\left[  \mathcal{A}\right]  $) and \newline$\left[
\mathcal{B}\right]  =%
\begin{cases}
1, & \text{if }\mathcal{B}\text{ is true};\\
0, & \text{if }\mathcal{B}\text{ is false}%
\end{cases}
$ (by the definition of $\left[  \mathcal{B}\right]  $).

But $\mathcal{A}$ and $\mathcal{B}$ are equivalent. Thus, $\mathcal{A}$ is
true (resp. false) if and only if $\mathcal{B}$ is true (resp. false). Hence,
$%
\begin{cases}
1, & \text{if }\mathcal{A}\text{ is true};\\
0, & \text{if }\mathcal{A}\text{ is false}%
\end{cases}
=%
\begin{cases}
1, & \text{if }\mathcal{B}\text{ is true};\\
0, & \text{if }\mathcal{B}\text{ is false}%
\end{cases}
$. Thus,%
\[
\left[  \mathcal{A}\right]  =%
\begin{cases}
1, & \text{if }\mathcal{A}\text{ is true};\\
0, & \text{if }\mathcal{A}\text{ is false}%
\end{cases}
=%
\begin{cases}
1, & \text{if }\mathcal{B}\text{ is true};\\
0, & \text{if }\mathcal{B}\text{ is false}%
\end{cases}
=\left[  \mathcal{B}\right]  .
\]
This solves Exercise \ref{exe.iverson-prop} \textbf{(a)}.

\textbf{(b)} Let $\mathcal{A}$ be any logical statement. Then, $\left(
\text{not }\mathcal{A}\right)  $ is true (resp. false) if and only if
$\mathcal{A}$ is false (resp. true). Hence,
\[%
\begin{cases}
1, & \text{if }\left(  \text{not }\mathcal{A}\right)  \text{ is true};\\
0, & \text{if }\left(  \text{not }\mathcal{A}\right)  \text{ is false}%
\end{cases}
=%
\begin{cases}
1, & \text{if }\mathcal{A}\text{ is false};\\
0, & \text{if }\mathcal{A}\text{ is true}%
\end{cases}
=%
\begin{cases}
0, & \text{if }\mathcal{A}\text{ is true};\\
1, & \text{if }\mathcal{A}\text{ is false}%
\end{cases}
.
\]
Now, the definition of $\left[  \text{not }\mathcal{A}\right]  $ shows that%
\[
\left[  \text{not }\mathcal{A}\right]  =%
\begin{cases}
1, & \text{if }\left(  \text{not }\mathcal{A}\right)  \text{ is true};\\
0, & \text{if }\left(  \text{not }\mathcal{A}\right)  \text{ is false}%
\end{cases}
=%
\begin{cases}
0, & \text{if }\mathcal{A}\text{ is true};\\
1, & \text{if }\mathcal{A}\text{ is false}%
\end{cases}
.
\]
Adding this equality to%
\[
\left[  \mathcal{A}\right]  =%
\begin{cases}
1, & \text{if }\mathcal{A}\text{ is true};\\
0, & \text{if }\mathcal{A}\text{ is false}%
\end{cases}
,
\]
we obtain%
\begin{align*}
\left[  \mathcal{A}\right]  +\left[  \text{not }\mathcal{A}\right]   &  =%
\begin{cases}
1, & \text{if }\mathcal{A}\text{ is true};\\
0, & \text{if }\mathcal{A}\text{ is false}%
\end{cases}
+%
\begin{cases}
0, & \text{if }\mathcal{A}\text{ is true};\\
1, & \text{if }\mathcal{A}\text{ is false}%
\end{cases}
=%
\begin{cases}
1+0, & \text{if }\mathcal{A}\text{ is true};\\
0+1, & \text{if }\mathcal{A}\text{ is false}%
\end{cases}
\\
&  =%
\begin{cases}
1, & \text{if }\mathcal{A}\text{ is true};\\
1, & \text{if }\mathcal{A}\text{ is false}%
\end{cases}
=1.
\end{align*}
Thus, $\left[  \text{not }\mathcal{A}\right]  =1-\left[  \mathcal{A}\right]
$. This solves Exercise \ref{exe.iverson-prop} \textbf{(b)}.

\textbf{(c)} Let $\mathcal{A}$ and $\mathcal{B}$ be two logical statements. We
must be in one of the following two cases:

\textit{Case 1:} The statement $\mathcal{A}$ is true.

\textit{Case 2:} The statement $\mathcal{A}$ is false.

Let us consider Case 1 first. In this case, the statement $\mathcal{A}$ is
true. Hence, the statement $\mathcal{A}\wedge\mathcal{B}$ is equivalent to the
statement $\mathcal{B}$. Thus, Exercise \ref{exe.iverson-prop} \textbf{(a)}
(applied to $\mathcal{A}\wedge\mathcal{B}$ instead of $\mathcal{A}$) shows
that $\left[  \mathcal{A}\wedge\mathcal{B}\right]  =\left[  \mathcal{B}%
\right]  $. But the definition of $\left[  \mathcal{A}\right]  $ yields
$\left[  \mathcal{A}\right]  =%
\begin{cases}
1, & \text{if }\mathcal{A}\text{ is true};\\
0, & \text{if }\mathcal{A}\text{ is false}%
\end{cases}
=1$ (since $\mathcal{A}$ is true). Hence, $\underbrace{\left[  \mathcal{A}%
\right]  }_{=1}\left[  \mathcal{B}\right]  =\left[  \mathcal{B}\right]  $.
Comparing this with $\left[  \mathcal{A}\wedge\mathcal{B}\right]  =\left[
\mathcal{B}\right]  $, we obtain $\left[  \mathcal{A}\wedge\mathcal{B}\right]
=\left[  \mathcal{A}\right]  \left[  \mathcal{B}\right]  $. Thus, Exercise
\ref{exe.iverson-prop} \textbf{(c)} is solved in Case 1.

Let us now consider Case 2. In this case, the statement $\mathcal{A}$ is
false. Hence, the statement $\mathcal{A}\wedge\mathcal{B}$ is false as well.
Thus, the definition of $\left[  \mathcal{A}\wedge\mathcal{B}\right]  $ yields
$\left[  \mathcal{A}\wedge\mathcal{B}\right]  =%
\begin{cases}
1, & \text{if }\mathcal{A}\wedge\mathcal{B}\text{ is true};\\
0, & \text{if }\mathcal{A}\wedge\mathcal{B}\text{ is false}%
\end{cases}
=0$ (since $\mathcal{A}\wedge\mathcal{B}$ is false). But the definition of
$\left[  \mathcal{A}\right]  $ yields $\left[  \mathcal{A}\right]  =%
\begin{cases}
1, & \text{if }\mathcal{A}\text{ is true};\\
0, & \text{if }\mathcal{A}\text{ is false}%
\end{cases}
=0$ (since $\mathcal{A}$ is false). Hence, $\underbrace{\left[  \mathcal{A}%
\right]  }_{=0}\left[  \mathcal{B}\right]  =0\left[  \mathcal{B}\right]  =0$.
Comparing this with $\left[  \mathcal{A}\wedge\mathcal{B}\right]  =0$, we
obtain $\left[  \mathcal{A}\wedge\mathcal{B}\right]  =\left[  \mathcal{A}%
\right]  \left[  \mathcal{B}\right]  $. Thus, Exercise \ref{exe.iverson-prop}
\textbf{(c)} is solved in Case 2.

We thus have solved Exercise \ref{exe.iverson-prop} \textbf{(c)} in both Cases
1 and 2. Hence, Exercise \ref{exe.iverson-prop} \textbf{(c)} always holds.

[\textit{Remark:} It is, of course, also possible to get a completely
straightforward solution to Exercise \ref{exe.iverson-prop} \textbf{(c)} by
distinguishing four cases, depending on which of the statements $\mathcal{A}$
and $\mathcal{B}$ are true.]

\textbf{(d)} It is easy to solve Exercise \ref{exe.iverson-prop} \textbf{(d)}
by a case distinction similarly to Exercise \ref{exe.iverson-prop}
\textbf{(c)}. However, since we have already solved parts \textbf{(b)} and
\textbf{(c)}, we can give a simpler solution:

Let $\mathcal{A}$ and $\mathcal{B}$ be two logical statements. One of de
Morgan's laws says that the statement $\left(  \text{not }\left(
\mathcal{A}\vee\mathcal{B}\right)  \right)  $ is equivalent to $\left(
\text{not }\mathcal{A}\right)  \wedge\left(  \text{not }\mathcal{B}\right)  $.
Hence, Exercise \ref{exe.iverson-prop} \textbf{(a)} (applied to $\left(
\text{not }\left(  \mathcal{A}\vee\mathcal{B}\right)  \right)  $ and $\left(
\text{not }\mathcal{A}\right)  \wedge\left(  \text{not }\mathcal{B}\right)  $
instead of $\mathcal{A}$ and $\mathcal{B}$) shows that%
\begin{align*}
\left[  \text{not }\left(  \mathcal{A}\vee\mathcal{B}\right)  \right]   &
=\left[  \left(  \text{not }\mathcal{A}\right)  \wedge\left(  \text{not
}\mathcal{B}\right)  \right] \\
&  =\underbrace{\left[  \text{not }\mathcal{A}\right]  }_{\substack{=1-\left[
\mathcal{A}\right]  \\\text{(by Exercise \ref{exe.iverson-prop} \textbf{(b)}%
)}}}\underbrace{\left[  \text{not }\mathcal{B}\right]  }_{\substack{=1-\left[
\mathcal{B}\right]  \\\text{(by Exercise \ref{exe.iverson-prop} \textbf{(b)}%
,}\\\text{applied to }\mathcal{B}\text{ instead of }\mathcal{A}\text{)}}}\\
&  \ \ \ \ \ \ \ \ \ \ \left(
\begin{array}
[c]{c}%
\text{by Exercise \ref{exe.iverson-prop} \textbf{(c)}, applied to}\\
\left(  \text{not }\mathcal{A}\right)  \text{ and }\left(  \text{not
}\mathcal{B}\right)  \text{ instead of }\mathcal{A}\text{ and }\mathcal{B}%
\end{array}
\right) \\
&  =\left(  1-\left[  \mathcal{A}\right]  \right)  \left(  1-\left[
\mathcal{B}\right]  \right)  =1-\left[  \mathcal{A}\right]  -\left[
\mathcal{B}\right]  +\left[  \mathcal{A}\right]  \left[  \mathcal{B}\right]  .
\end{align*}
But Exercise \ref{exe.iverson-prop} \textbf{(b)} (applied to $\mathcal{A}%
\vee\mathcal{B}$ instead of $\mathcal{A}$) shows that $\left[  \text{not
}\left(  \mathcal{A}\vee\mathcal{B}\right)  \right]  =1-\left[  \mathcal{A}%
\vee\mathcal{B}\right]  $. Hence,%
\[
\left[  \mathcal{A}\vee\mathcal{B}\right]  =1-\underbrace{\left[  \text{not
}\left(  \mathcal{A}\vee\mathcal{B}\right)  \right]  }_{=1-\left[
\mathcal{A}\right]  -\left[  \mathcal{B}\right]  +\left[  \mathcal{A}\right]
\left[  \mathcal{B}\right]  }=1-\left(  1-\left[  \mathcal{A}\right]  -\left[
\mathcal{B}\right]  +\left[  \mathcal{A}\right]  \left[  \mathcal{B}\right]
\right)  =\left[  \mathcal{A}\right]  +\left[  \mathcal{B}\right]  -\left[
\mathcal{A}\right]  \left[  \mathcal{B}\right]  .
\]
This solves Exercise \ref{exe.iverson-prop} \textbf{(d)}.

\textbf{(e)} Let $\mathcal{A}$, $\mathcal{B}$ and $\mathcal{C}$ be three
logical statements. Then, Exercise \ref{exe.iverson-prop} \textbf{(d)}
(applied to $\mathcal{A}\vee\mathcal{B}$ and $\mathcal{C}$ instead of
$\mathcal{A}$ and $\mathcal{B}$) shows that%
\begin{align*}
\left[  \left(  \mathcal{A}\vee\mathcal{B}\right)  \vee\mathcal{C}\right]   &
=\underbrace{\left[  \mathcal{A}\vee\mathcal{B}\right]  }_{\substack{=\left[
\mathcal{A}\right]  +\left[  \mathcal{B}\right]  -\left[  \mathcal{A}\right]
\left[  \mathcal{B}\right]  \\\text{(by Exercise \ref{exe.iverson-prop}
\textbf{(d)})}}}+\left[  \mathcal{C}\right]  -\underbrace{\left[
\mathcal{A}\vee\mathcal{B}\right]  }_{\substack{=\left[  \mathcal{A}\right]
+\left[  \mathcal{B}\right]  -\left[  \mathcal{A}\right]  \left[
\mathcal{B}\right]  \\\text{(by Exercise \ref{exe.iverson-prop} \textbf{(d)}%
)}}}\left[  \mathcal{C}\right] \\
&  =\left(  \left[  \mathcal{A}\right]  +\left[  \mathcal{B}\right]  -\left[
\mathcal{A}\right]  \left[  \mathcal{B}\right]  \right)  +\left[
\mathcal{C}\right]  -\left(  \left[  \mathcal{A}\right]  +\left[
\mathcal{B}\right]  -\left[  \mathcal{A}\right]  \left[  \mathcal{B}\right]
\right)  \left[  \mathcal{C}\right] \\
&  =\left[  \mathcal{A}\right]  +\left[  \mathcal{B}\right]  +\left[
\mathcal{C}\right]  -\left[  \mathcal{A}\right]  \left[  \mathcal{B}\right]
-\left[  \mathcal{A}\right]  \left[  \mathcal{C}\right]  -\left[
\mathcal{B}\right]  \left[  \mathcal{C}\right]  +\left[  \mathcal{A}\right]
\left[  \mathcal{B}\right]  \left[  \mathcal{C}\right]  .
\end{align*}
But the statement $\mathcal{A}\vee\mathcal{B}\vee\mathcal{C}$ is equivalent to
$\left(  \mathcal{A}\vee\mathcal{B}\right)  \vee\mathcal{C}$. Hence, Exercise
\ref{exe.iverson-prop} \textbf{(a)} (applied to $\mathcal{A}\vee
\mathcal{B}\vee\mathcal{C}$ and $\left(  \mathcal{A}\vee\mathcal{B}\right)
\vee\mathcal{C}$ instead of $\mathcal{A}$ and $\mathcal{B}$) shows that%
\begin{align*}
\left[  \mathcal{A}\vee\mathcal{B}\vee\mathcal{C}\right]   &  =\left[  \left(
\mathcal{A}\vee\mathcal{B}\right)  \vee\mathcal{C}\right] \\
&  =\left[  \mathcal{A}\right]  +\left[  \mathcal{B}\right]  +\left[
\mathcal{C}\right]  -\left[  \mathcal{A}\right]  \left[  \mathcal{B}\right]
-\left[  \mathcal{A}\right]  \left[  \mathcal{C}\right]  -\left[
\mathcal{B}\right]  \left[  \mathcal{C}\right]  +\left[  \mathcal{A}\right]
\left[  \mathcal{B}\right]  \left[  \mathcal{C}\right]  .
\end{align*}
This solves Exercise \ref{exe.iverson-prop} \textbf{(e)}.
\end{proof}

\subsection{Solution to Exercise \ref{exe.pie.specialize}}

\begin{vershort}
\begin{proof}
[Proof of Theorem \ref{thm.pie.bonfer}.]Clearly, $m\geq0$. Hence, Theorem
\ref{thm.pie.merged} (applied to $k=0$) yields%
\begin{align*}
\left(  -1\right)  ^{m}\sum_{s\in S}\dbinom{c\left(  s\right)  }{0}%
\dbinom{c\left(  s\right)  -0-1}{m-0}  &  =\sum_{\substack{I\subseteq
G;\\\left\vert I\right\vert \leq m}}\left(  -1\right)  ^{\left\vert
I\right\vert }\underbrace{\dbinom{\left\vert I\right\vert }{0}}%
_{\substack{=1\\\text{(by Proposition \ref{prop.binom.00} \textbf{(a)}%
}\\\text{(applied to }\left\vert I\right\vert \text{ instead of }m\text{))}%
}}\left\vert \bigcap_{i\in I}A_{i}\right\vert \\
&  =\sum_{\substack{I\subseteq G;\\\left\vert I\right\vert \leq m}}\left(
-1\right)  ^{\left\vert I\right\vert }\left\vert \bigcap_{i\in I}%
A_{i}\right\vert .
\end{align*}
Hence,%
\begin{align*}
\sum_{\substack{I\subseteq G;\\\left\vert I\right\vert \leq m}}\left(
-1\right)  ^{\left\vert I\right\vert }\left\vert \bigcap_{i\in I}%
A_{i}\right\vert  &  =\left(  -1\right)  ^{m}\sum_{s\in S}\underbrace{\dbinom
{c\left(  s\right)  }{0}}_{\substack{=1\\\text{(by Proposition
\ref{prop.binom.00} \textbf{(a)}}\\\text{(applied to }c\left(  s\right)
\text{ instead of }m\text{))}}}\underbrace{\dbinom{c\left(  s\right)
-0-1}{m-0}}_{=\dbinom{c\left(  s\right)  -1}{m}}\\
&  =\left(  -1\right)  ^{m}\sum_{s\in S}\dbinom{c\left(  s\right)  -1}{m}.
\end{align*}
This proves Theorem \ref{thm.pie.bonfer}.
\end{proof}
\end{vershort}

\begin{verlong}
\begin{proof}
[Proof of Theorem \ref{thm.pie.bonfer}.]Clearly, $m\geq0$. Hence, Theorem
\ref{thm.pie.merged} (applied to $k=0$) yields%
\begin{align*}
\left(  -1\right)  ^{m}\sum_{s\in S}\dbinom{c\left(  s\right)  }{0}%
\dbinom{c\left(  s\right)  -0-1}{m-0}  &  =\sum_{\substack{I\subseteq
G;\\\left\vert I\right\vert \leq m}}\left(  -1\right)  ^{\left\vert
I\right\vert }\underbrace{\dbinom{\left\vert I\right\vert }{0}}%
_{\substack{=1\\\text{(by Proposition \ref{prop.binom.00} \textbf{(a)}%
}\\\text{(applied to }\left\vert I\right\vert \text{ instead of }m\text{))}%
}}\left\vert \bigcap_{i\in I}A_{i}\right\vert \\
&  =\sum_{\substack{I\subseteq G;\\\left\vert I\right\vert \leq m}}\left(
-1\right)  ^{\left\vert I\right\vert }\left\vert \bigcap_{i\in I}%
A_{i}\right\vert .
\end{align*}
Hence,%
\begin{align*}
\sum_{\substack{I\subseteq G;\\\left\vert I\right\vert \leq m}}\left(
-1\right)  ^{\left\vert I\right\vert }\left\vert \bigcap_{i\in I}%
A_{i}\right\vert  &  =\left(  -1\right)  ^{m}\sum_{s\in S}\underbrace{\dbinom
{c\left(  s\right)  }{0}}_{\substack{=1\\\text{(by Proposition
\ref{prop.binom.00} \textbf{(a)}}\\\text{(applied to }c\left(  s\right)
\text{ instead of }m\text{))}}}\dbinom{c\left(  s\right)  -0-1}{m-0}\\
&  =\left(  -1\right)  ^{m}\sum_{s\in S}\underbrace{\dbinom{c\left(  s\right)
-0-1}{m-0}}_{\substack{=\dbinom{c\left(  s\right)  -1}{m}\\\text{(since
}c\left(  s\right)  -0=c\left(  s\right)  \text{ and }m-0=m\text{)}}}=\left(
-1\right)  ^{m}\sum_{s\in S}\dbinom{c\left(  s\right)  -1}{m}.
\end{align*}
This proves Theorem \ref{thm.pie.bonfer}.
\end{proof}
\end{verlong}

\begin{proof}
[Proof of Theorem \ref{thm.pie.jordan}.]Clearly, $\left\vert G\right\vert
\in\mathbb{N}$ (since the set $G$ is finite). Define an $m\in\mathbb{N}$ by
$m=\max\left\{  \left\vert G\right\vert ,k\right\}  $. Thus, $m=\max\left\{
\left\vert G\right\vert ,k\right\}  \geq\left\vert G\right\vert $ and
$m=\max\left\{  \left\vert G\right\vert ,k\right\}  \geq k$. Also, every
subset $I$ of $G$ satisfies $\left\vert I\right\vert \leq m$%
\ \ \ \ \footnote{\textit{Proof.} Let $I$ be a subset of $G$. Thus,
$\left\vert I\right\vert \leq\left\vert G\right\vert \leq m$ (since
$m\geq\left\vert G\right\vert $). Qed.}. Hence, we have the following equality
of summation signs:%
\begin{equation}
\sum_{\substack{I\subseteq G;\\\left\vert I\right\vert \leq m}}=\sum
_{I\subseteq G}. \label{pf.thm.pie.jordan.sum=sum}%
\end{equation}

For each $s\in S$, let $c\left(  s\right)  $ denote the number of $i\in G$
satisfying $s\in A_{i}$.

Recall that
\[
S_{k}=\left\{  s\in S\ \mid\ \text{the number of }i\in G\text{ satisfying
}s\in A_{i}\text{ equals }k\right\}  .
\]
Hence, for each $s\in S$, we have the following chain of logical equivalences:%
\begin{align}
\left(  s\in S_{k}\right)  \  &  \Longleftrightarrow\ \left(
\underbrace{\text{the number of }i\in G\text{ satisfying }s\in A_{i}%
}_{\substack{=c\left(  s\right)  \\\text{(since }c\left(  s\right)  \text{ is
the number of }i\in G\text{ satisfying }s\in A_{i}\\\text{(by the definition
of }c\left(  s\right)  \text{))}}}\text{ equals }k\right) \nonumber\\
\  &  \Longleftrightarrow\ \left(  c\left(  s\right)  \text{ equals }k\right)
\nonumber\\
\  &  \Longleftrightarrow\ \left(  c\left(  s\right)  =k\right)  .
\label{pf.thm.pie.jordan.equiv}%
\end{align}

Observe that $m-k\in\mathbb{N}$ (since $m\geq k$).

Next, we notice the following:

\begin{statement}
\textit{Observation 1:} For any $s\in S$, we have
\[
\dbinom{c\left(  s\right)  }{k}\dbinom{c\left(  s\right)  -k-1}{m-k}=\left(
-1\right)  ^{m-k}\left[  s\in S_{k}\right]  .
\]

\end{statement}

[\textit{Proof of Observation 1:} Let $s\in S$. Clearly, $c\left(  s\right)
\in\mathbb{N}$. We are in one of the following three cases:

\textit{Case 1:} We have $c\left(  s\right)  <k$.

\textit{Case 2:} We have $c\left(  s\right)  =k$.

\textit{Case 3:} We have $c\left(  s\right)  >k$.

Let us first consider Case 1. In this case, we have $c\left(  s\right)  <k$.
Thus, Proposition \ref{prop.binom.0} (applied to $c\left(  s\right)  $ and $k$
instead of $m$ and $n$) yields $\dbinom{c\left(  s\right)  }{k}=0$. Hence,
$\underbrace{\dbinom{c\left(  s\right)  }{k}}_{=0}\dbinom{c\left(  s\right)
-k-1}{m-k}=0$.

But $c\left(  s\right)  \neq k$ (since $c\left(  s\right)  <k$). Thus, we
don't have $c\left(  s\right)  =k$. Hence, we don't have $s\in S_{k}$ (by the
equivalence (\ref{pf.thm.pie.jordan.equiv})). Thus, $\left[  s\in
S_{k}\right]  =0$. Hence, $\left(  -1\right)  ^{m-k}\underbrace{\left[  s\in
S_{k}\right]  }_{=0}=0$. Comparing this with $\dbinom{c\left(  s\right)  }%
{k}\dbinom{c\left(  s\right)  -k-1}{m-k}=0$, we obtain $\dbinom{c\left(
s\right)  }{k}\dbinom{c\left(  s\right)  -k-1}{m-k}=\left(  -1\right)
^{m-k}\left[  s\in S_{k}\right]  $. Hence, Observation 1 is proven in Case 1.

\begin{vershort}
Let us now consider Case 2. In this case, we have $c\left(  s\right)  =k$.
Hence,%
\[
\dbinom{c\left(  s\right)  }{k}\dbinom{c\left(  s\right)  -k-1}{m-k}%
=\underbrace{\dbinom{k}{k}}_{\substack{=1\\\text{(by Proposition
\ref{prop.binom.mm}}\\\text{(applied to }k\text{ instead of }m\text{))}%
}}\underbrace{\dbinom{k-k-1}{m-k}}_{=\dbinom{-1}{m-k}}=\dbinom{-1}%
{m-k}=\left(  -1\right)  ^{m-k}%
\]
(by Corollary \ref{cor.binom.-1} (applied to $n=m-k$)).
\end{vershort}

\begin{verlong}
Let us now consider Case 2. In this case, we have $c\left(  s\right)  =k$.
Hence,%
\begin{align*}
\dbinom{c\left(  s\right)  }{k}\dbinom{c\left(  s\right)  -k-1}{m-k}  &
=\underbrace{\dbinom{k}{k}}_{\substack{=1\\\text{(by Proposition
\ref{prop.binom.mm}}\\\text{(applied to }k\text{ instead of }m\text{))}%
}}\underbrace{\dbinom{k-k-1}{m-k}}_{\substack{=\dbinom{-1}{m-k}\\\text{(since
}k-k-1=-1\text{)}}}\\
&  =\dbinom{-1}{m-k}=\left(  -1\right)  ^{m-k}%
\end{align*}
(by Corollary \ref{cor.binom.-1} (applied to $n=m-k$)).
\end{verlong}

But $c\left(  s\right)  =k$. Hence, $s\in S_{k}$ (by the equivalence
(\ref{pf.thm.pie.jordan.equiv})). Thus, $\left[  s\in S_{k}\right]  =1$.
Hence, $\left(  -1\right)  ^{m-k}\underbrace{\left[  s\in S_{k}\right]  }%
_{=1}=\left(  -1\right)  ^{m-k}$. Comparing this with $\dbinom{c\left(
s\right)  }{k}\dbinom{c\left(  s\right)  -k-1}{m-k}=\left(  -1\right)  ^{m-k}%
$, we obtain $\dbinom{c\left(  s\right)  }{k}\dbinom{c\left(  s\right)
-k-1}{m-k}=\left(  -1\right)  ^{m-k}\left[  s\in S_{k}\right]  $. Hence,
Observation 1 is proven in Case 2.

Let us finally consider Case 3. In this case, we have $c\left(  s\right)  >k$.
Thus, $c\left(  s\right)  \geq k+1$ (since $c\left(  s\right)  $ and $k$ are
integers), so that $c\left(  s\right)  -k\geq1$ and thus $c\left(  s\right)
-k-1\in\mathbb{N}$. But $c\left(  s\right)  $ is the number of $i\in G$
satisfying $s\in A_{i}$. Hence,
\[
c\left(  s\right)  =\left(  \text{the number of }i\in G\text{ satisfying }s\in
A_{i}\right)  =\left\vert \underbrace{\left\{  i\in G\ \mid\ s\in
A_{i}\right\}  }_{\subseteq G}\right\vert \leq\left\vert G\right\vert \leq m
\]
(since $m\geq\left\vert G\right\vert $). Hence, $\underbrace{c\left(
s\right)  }_{\leq m}-k-1\leq m-k-1<m-k$. Therefore, Proposition
\ref{prop.binom.0} (applied to $c\left(  s\right)  -k-1$ and $m-k$ instead of
$m$ and $n$) yields $\dbinom{c\left(  s\right)  -k-1}{m-k}=0$. Hence,
$\dbinom{c\left(  s\right)  }{k}\underbrace{\dbinom{c\left(  s\right)
-k-1}{m-k}}_{=0}=0$.

But $c\left(  s\right)  \neq k$ (since $c\left(  s\right)  >k$). Thus, we
don't have $c\left(  s\right)  =k$. Hence, we don't have $s\in S_{k}$ (by the
equivalence (\ref{pf.thm.pie.jordan.equiv})). Thus, $\left[  s\in
S_{k}\right]  =0$. Hence, $\left(  -1\right)  ^{m-k}\underbrace{\left[  s\in
S_{k}\right]  }_{=0}=0$. Comparing this with $\dbinom{c\left(  s\right)  }%
{k}\dbinom{c\left(  s\right)  -k-1}{m-k}=0$, we obtain $\dbinom{c\left(
s\right)  }{k}\dbinom{c\left(  s\right)  -k-1}{m-k}=\left(  -1\right)
^{m-k}\left[  s\in S_{k}\right]  $. Hence, Observation 1 is proven in Case 3.

We have now proven Observation 1 in each of the three Cases 1, 2 and 3. Hence,
Observation 1 always holds.]

Now, Theorem \ref{thm.pie.merged} yields%
\begin{align*}
\left(  -1\right)  ^{m}\sum_{s\in S}\dbinom{c\left(  s\right)  }{k}%
\dbinom{c\left(  s\right)  -k-1}{m-k}  &  =\underbrace{\sum
_{\substack{I\subseteq G;\\\left\vert I\right\vert \leq m}}}_{\substack{=\sum
_{I\subseteq G}\\\text{(by (\ref{pf.thm.pie.jordan.sum=sum}))}}}\left(
-1\right)  ^{\left\vert I\right\vert }\dbinom{\left\vert I\right\vert }%
{k}\left\vert \bigcap_{i\in I}A_{i}\right\vert \\
&  =\sum_{I\subseteq G}\left(  -1\right)  ^{\left\vert I\right\vert }%
\dbinom{\left\vert I\right\vert }{k}\left\vert \bigcap_{i\in I}A_{i}%
\right\vert .
\end{align*}
Thus,%
\begin{align}
&  \sum_{I\subseteq G}\left(  -1\right)  ^{\left\vert I\right\vert }%
\dbinom{\left\vert I\right\vert }{k}\left\vert \bigcap_{i\in I}A_{i}%
\right\vert \nonumber\\
&  =\left(  -1\right)  ^{m}\sum_{s\in S}\underbrace{\dbinom{c\left(  s\right)
}{k}\dbinom{c\left(  s\right)  -k-1}{m-k}}_{\substack{=\left(  -1\right)
^{m-k}\left[  s\in S_{k}\right]  \\\text{(by Observation 1)}}}\nonumber\\
&  =\left(  -1\right)  ^{m}\sum_{s\in S}\left(  -1\right)  ^{m-k}\left[  s\in
S_{k}\right]  =\underbrace{\left(  -1\right)  ^{m}\left(  -1\right)  ^{m-k}%
}_{\substack{=\left(  -1\right)  ^{m+\left(  m-k\right)  }=\left(  -1\right)
^{k}\\\text{(since }m+\left(  m-k\right)  =2m-k\equiv-k\equiv
k\operatorname{mod}2\text{)}}}\sum_{s\in S}\left[  s\in S_{k}\right]
\nonumber\\
&  =\left(  -1\right)  ^{k}\sum_{s\in S}\left[  s\in S_{k}\right]  .
\label{pf.thm.pie.jordan.almost}%
\end{align}
Multiplying both sides of this equality by $\left(  -1\right)  ^{k}$, we find%
\begin{align}
&  \left(  -1\right)  ^{k}\sum_{I\subseteq G}\left(  -1\right)  ^{\left\vert
I\right\vert }\dbinom{\left\vert I\right\vert }{k}\left\vert \bigcap_{i\in
I}A_{i}\right\vert \nonumber\\
&  =\underbrace{\left(  -1\right)  ^{k}\left(  -1\right)  ^{k}}%
_{\substack{=\left(  \left(  -1\right)  \left(  -1\right)  \right)  ^{k}%
=1^{k}\\\text{(since }\left(  -1\right)  \left(  -1\right)  =1\text{)}}%
}\sum_{s\in S}\left[  s\in S_{k}\right]  =\underbrace{1^{k}}_{=1}\sum_{s\in
S}\left[  s\in S_{k}\right] \nonumber\\
&  =\sum_{s\in S}\left[  s\in S_{k}\right]  . \label{pf.thm.pie.jordan.almos2}%
\end{align}

But $S_{k}$ is a subset of $S$ (by the definition of $S_{k}$). Hence, Lemma
\ref{lem.iverson.card} (applied to $T=S_{k}$) yields%
\begin{align*}
\left\vert S_{k}\right\vert  &  =\sum_{s\in S}\left[  s\in S_{k}\right]
=\left(  -1\right)  ^{k}\sum_{I\subseteq G}\left(  -1\right)  ^{\left\vert
I\right\vert }\dbinom{\left\vert I\right\vert }{k}\left\vert \bigcap_{i\in
I}A_{i}\right\vert \ \ \ \ \ \ \ \ \ \ \left(  \text{by
(\ref{pf.thm.pie.jordan.almos2})}\right) \\
&  =\sum_{I\subseteq G}\underbrace{\left(  -1\right)  ^{\left\vert
I\right\vert }\left(  -1\right)  ^{k}}_{\substack{=\left(  -1\right)
^{\left\vert I\right\vert +k}=\left(  -1\right)  ^{\left\vert I\right\vert
-k}\\\text{(since }\left\vert I\right\vert +k\equiv\left\vert I\right\vert
-k\operatorname{mod}2\text{)}}}\dbinom{\left\vert I\right\vert }{k}\left\vert
\bigcap_{i\in I}A_{i}\right\vert =\sum_{I\subseteq G}\left(  -1\right)
^{\left\vert I\right\vert -k}\dbinom{\left\vert I\right\vert }{k}\left\vert
\bigcap_{i\in I}A_{i}\right\vert .
\end{align*}
This proves Theorem \ref{thm.pie.jordan}.
\end{proof}

\begin{proof}
[Proof of Theorem \ref{thm.pie.nonunion}.]We have $\bigcup_{i\in G}%
A_{i}\subseteq S$ (since all the $A_{i}$ are subsets of $S$). Hence,%
\begin{align}
\bigcup_{i\in G}A_{i}  &  =\underbrace{\left(  \bigcup_{i\in G}A_{i}\right)
}_{\substack{=\left\{  s\ \mid\ \text{there exists an }i\in G\text{ such that
}s\in A_{i}\right\}  \\\text{(by the definition of the union }\bigcup_{i\in
G}A_{i}\text{)}}}\cap S\nonumber\\
&  =\left\{  s\ \mid\ \text{there exists an }i\in G\text{ such that }s\in
A_{i}\right\}  \cap S\nonumber\\
&  =\left\{  s\in S\ \mid\ \text{there exists an }i\in G\text{ such that }s\in
A_{i}\right\}  . \label{pf.thm.pie.nonunion.union}%
\end{align}

Define a $k\in\mathbb{N}$ by $k=0$. Define the set $S_{k}$ as in Theorem
\ref{thm.pie.jordan}. Then,
\begin{align*}
S_{k}  &  =\left\{  s\in S\ \mid\ \text{the number of }i\in G\text{ satisfying
}s\in A_{i}\text{ equals }k\right\} \\
&  =\left\{  s\in S\ \mid\ \underbrace{\text{the number of }i\in G\text{
satisfying }s\in A_{i}\text{ equals }0}_{\Longleftrightarrow\ \left(
\text{there exists no }i\in G\text{ such that }s\in A_{i}\right)  }\right\} \\
&  \ \ \ \ \ \ \ \ \ \ \left(  \text{since }k=0\right) \\
&  =\left\{  s\in S\ \mid\ \text{there exists no }i\in G\text{ such that }s\in
A_{i}\right\} \\
&  =S\setminus\underbrace{\left\{  s\in S\ \mid\ \text{there exists an }i\in
G\text{ such that }s\in A_{i}\right\}  }_{\substack{=\bigcup_{i\in G}%
A_{i}\\\text{(by (\ref{pf.thm.pie.nonunion.union}))}}}\\
&  =S\setminus\left(  \bigcup_{i\in G}A_{i}\right)  .
\end{align*}
Hence, $\left\vert S_{k}\right\vert =\left\vert S\setminus\left(
\bigcup_{i\in G}A_{i}\right)  \right\vert $. Therefore,%
\begin{align*}
\left\vert S\setminus\left(  \bigcup_{i\in G}A_{i}\right)  \right\vert  &
=\left\vert S_{k}\right\vert =\sum_{I\subseteq G}\underbrace{\left(
-1\right)  ^{\left\vert I\right\vert -k}}_{\substack{=\left(  -1\right)
^{\left\vert I\right\vert -0}\\\text{(since }k=0\text{)}}}\underbrace{\dbinom
{\left\vert I\right\vert }{k}}_{\substack{=\dbinom{\left\vert I\right\vert
}{0}\\\text{(since }k=0\text{)}}}\left\vert \bigcap_{i\in I}A_{i}\right\vert
\\
&  \ \ \ \ \ \ \ \ \ \ \left(  \text{by Theorem \ref{thm.pie.jordan}}\right)
\\
&  =\sum_{I\subseteq G}\underbrace{\left(  -1\right)  ^{\left\vert
I\right\vert -0}}_{=\left(  -1\right)  ^{\left\vert I\right\vert }%
}\underbrace{\dbinom{\left\vert I\right\vert }{0}}_{\substack{=1\\\text{(by
Proposition \ref{prop.binom.00} \textbf{(a)}}\\\text{(applied to }\left\vert
I\right\vert \text{ instead of }m\text{))}}}\left\vert \bigcap_{i\in I}%
A_{i}\right\vert \\
&  =\sum_{I\subseteq G}\left(  -1\right)  ^{\left\vert I\right\vert
}\left\vert \bigcap_{i\in I}A_{i}\right\vert .
\end{align*}
This proves Theorem \ref{thm.pie.nonunion}.
\end{proof}

\begin{proof}
[Proof of Theorem \ref{thm.pie.union}.]Let $S$ denote the set $\bigcup_{i\in
G}A_{i}$. Then, $S$ is the union of finitely many finite sets (since the set
$G$ is finite, and since each of the sets $A_{i}$ is finite), and thus itself
is a finite set. Moreover, $S=\bigcup_{i\in G}A_{i}$; thus, $A_{i}$ is a
subset of $S$ for each $i\in G$.

We define the intersection $\bigcap_{i\in\varnothing}A_{i}$ (which would
otherwise be undefined, since $\varnothing$ is the empty set) to mean the set
$S$. (Thus, $\bigcap_{i\in I}A_{i}$ is defined for any subset $I$ of $G$, not
just for nonempty subsets $I$.)

We have $\bigcap_{i\in\varnothing}A_{i}=S$ (since we have defined
$\bigcap_{i\in\varnothing}A_{i}$ to be $S$).

Theorem \ref{thm.pie.nonunion} yields%
\begin{align*}
\left\vert S\setminus\left(  \bigcup_{i\in G}A_{i}\right)  \right\vert  &
=\sum_{I\subseteq G}\left(  -1\right)  ^{\left\vert I\right\vert }\left\vert
\bigcap_{i\in I}A_{i}\right\vert \\
&  =\sum_{\substack{I\subseteq G;\\I\neq\varnothing}}\left(  -1\right)
^{\left\vert I\right\vert }\left\vert \bigcap_{i\in I}A_{i}\right\vert
+\underbrace{\left(  -1\right)  ^{\left\vert \varnothing\right\vert }%
}_{\substack{=\left(  -1\right)  ^{0}\\\text{(since }\left\vert \varnothing
\right\vert =0\text{)}}}\left\vert \underbrace{\bigcap_{i\in\varnothing}A_{i}%
}_{=S}\right\vert \\
&  \ \ \ \ \ \ \ \ \ \ \left(
\begin{array}
[c]{c}%
\text{here, we have split off the addend for }I=\varnothing\text{ from the
sum,}\\
\text{since }\varnothing\text{ is a subset of }G
\end{array}
\right) \\
&  =\sum_{\substack{I\subseteq G;\\I\neq\varnothing}}\left(  -1\right)
^{\left\vert I\right\vert }\left\vert \bigcap_{i\in I}A_{i}\right\vert
+\underbrace{\left(  -1\right)  ^{0}}_{=1}\left\vert S\right\vert
=\sum_{\substack{I\subseteq G;\\I\neq\varnothing}}\left(  -1\right)
^{\left\vert I\right\vert }\left\vert \bigcap_{i\in I}A_{i}\right\vert
+\left\vert S\right\vert .
\end{align*}
Comparing this with%
\[
\left\vert S\setminus\underbrace{\left(  \bigcup_{i\in G}A_{i}\right)  }%
_{=S}\right\vert =\left\vert \underbrace{S\setminus S}_{=\varnothing
}\right\vert =\left\vert \varnothing\right\vert =0,
\]
we obtain%
\[
0=\sum_{\substack{I\subseteq G;\\I\neq\varnothing}}\left(  -1\right)
^{\left\vert I\right\vert }\left\vert \bigcap_{i\in I}A_{i}\right\vert
+\left\vert S\right\vert .
\]
Solving this equation for $\left\vert S\right\vert $, we find%
\[
\left\vert S\right\vert =-\sum_{\substack{I\subseteq G;\\I\neq\varnothing
}}\left(  -1\right)  ^{\left\vert I\right\vert }\left\vert \bigcap_{i\in
I}A_{i}\right\vert =\sum_{\substack{I\subseteq G;\\I\neq\varnothing
}}\underbrace{\left(  -\left(  -1\right)  ^{\left\vert I\right\vert }\right)
}_{=\left(  -1\right)  ^{\left\vert I\right\vert -1}}\left\vert \bigcap_{i\in
I}A_{i}\right\vert =\sum_{\substack{I\subseteq G;\\I\neq\varnothing}}\left(
-1\right)  ^{\left\vert I\right\vert -1}\left\vert \bigcap_{i\in I}%
A_{i}\right\vert .
\]
In view of $S=\bigcup_{i\in G}A_{i}$, this rewrites as
\[
\left\vert \bigcup_{i\in G}A_{i}\right\vert =\sum_{\substack{I\subseteq
G;\\I\neq\varnothing}}\left(  -1\right)  ^{\left\vert I\right\vert
-1}\left\vert \bigcap_{i\in I}A_{i}\right\vert .
\]
This proves Theorem \ref{thm.pie.union}.
\end{proof}

\begin{proof}
[Solution to Exercise \ref{exe.pie.specialize}.]We have proven Theorem
\ref{thm.pie.union}, Theorem \ref{thm.pie.nonunion}, Theorem
\ref{thm.pie.jordan} and Theorem \ref{thm.pie.bonfer}. Thus, Exercise
\ref{exe.pie.specialize} is solved.
\end{proof}

\subsection{Solution to Exercise \ref{exe.multichoose}}

Exercise \ref{exe.multichoose} is one of the most fundamental results in
combinatorics. A combinatorial proof is sketched, e.g., in \cite[proof of
Proposition 13.3]{Galvin} and in \cite[proof of Identity 143]{BenQui03}. We
shall give a different proof, using induction instead.

Let us first state a basic lemma about sets:

\begin{lemma}
\label{lem.prodrule.prod-assM}Let $M$ be a positive integer. For every
$i\in\left\{  1,2,\ldots,M\right\}  $, let $Z_{i}$ be a set. Then, the map
\begin{align*}
Z_{1}\times Z_{2}\times\cdots\times Z_{M}  &  \rightarrow\left(  Z_{1}\times
Z_{2}\times\cdots\times Z_{M-1}\right)  \times Z_{M},\\
\left(  s_{1},s_{2},\ldots,s_{M}\right)   &  \mapsto\left(  \left(
s_{1},s_{2},\ldots,s_{M-1}\right)  ,s_{M}\right)
\end{align*}
is a bijection.
\end{lemma}

\begin{vershort}
\begin{proof}
[Proof of Lemma \ref{lem.prodrule.prod-assM}.]This map is the canonical
bijection $Z_{1}\times Z_{2}\times\cdots\times Z_{M}\rightarrow\left(
Z_{1}\times Z_{2}\times\cdots\times Z_{M-1}\right)  \times Z_{M}$. Its inverse
map sends each $\left(  \left(  s_{1},s_{2},\ldots,s_{M-1}\right)  ,t\right)
\in\left(  Z_{1}\times Z_{2}\times\cdots\times Z_{M-1}\right)  \times Z_{M}$
to $\left(  s_{1},s_{2},\ldots,s_{M-1},t\right)  \in Z_{1}\times Z_{2}%
\times\cdots\times Z_{M}$.
\end{proof}
\end{vershort}

\begin{verlong}
\begin{proof}
[Proof of Lemma \ref{lem.prodrule.prod-assM}.]We define a map%
\[
\Phi:\left(  Z_{1}\times Z_{2}\times\cdots\times Z_{M-1}\right)  \times
Z_{M}\rightarrow Z_{1}\times Z_{2}\times\cdots\times Z_{M}%
\]
by%
\begin{equation}
\left(
\begin{array}
[c]{l}%
\Phi\left(  \left(  s_{1},s_{2},\ldots,s_{M-1}\right)  ,t\right)  =\left(
s_{1},s_{2},\ldots,s_{M-1},t\right) \\
\ \ \ \ \ \ \ \ \ \ \text{for every }\left(  \left(  s_{1},s_{2}%
,\ldots,s_{M-1}\right)  ,t\right)  \in\left(  Z_{1}\times Z_{2}\times
\cdots\times Z_{M-1}\right)  \times Z_{M}%
\end{array}
\right)  . \label{pf.lem.prodrule.S.Phidef}%
\end{equation}

We also define a map%
\[
\Psi:Z_{1}\times Z_{2}\times\cdots\times Z_{M}\rightarrow\left(  Z_{1}\times
Z_{2}\times\cdots\times Z_{M-1}\right)  \times Z_{M}%
\]
by%
\begin{equation}
\left(
\begin{array}
[c]{l}%
\Psi\left(  s_{1},s_{2},\ldots,s_{M}\right)  =\left(  \left(  s_{1}%
,s_{2},\ldots,s_{M-1}\right)  ,s_{M}\right) \\
\ \ \ \ \ \ \ \ \ \ \text{for every }\left(  s_{1},s_{2},\ldots,s_{M}\right)
\in Z_{1}\times Z_{2}\times\cdots\times Z_{M}%
\end{array}
\right)  . \label{pf.lem.prodrule.S.Psidef}%
\end{equation}

We now claim that
\begin{equation}
\Psi\left(  s_{1},s_{2},\ldots,s_{M-1},t\right)  =\left(  \left(  s_{1}%
,s_{2},\ldots,s_{M-1}\right)  ,t\right)  \label{pf.lem.prodrule.S.Psi}%
\end{equation}
for every $\left(  s_{1},s_{2},\ldots,s_{M-1}\right)  \in Z_{1}\times
Z_{2}\times\cdots\times Z_{M-1}$ and $t\in Z_{M}$.

[\textit{Proof of (\ref{pf.lem.prodrule.S.Psi}):} Let $\left(  s_{1}%
,s_{2},\ldots,s_{M-1}\right)  \in Z_{1}\times Z_{2}\times\cdots\times Z_{M-1}$
and $t\in Z_{M}$. We have $\left(  s_{1},s_{2},\ldots,s_{M-1}\right)  \in
Z_{1}\times Z_{2}\times\cdots\times Z_{M-1}$; in other words,
\begin{equation}
s_{i}\in Z_{i}\ \ \ \ \ \ \ \ \ \ \text{for every }i\in\left\{  1,2,\ldots
,M-1\right\}  . \label{pf.lem.prodrule.S.Psi.pf.1}%
\end{equation}

For every $i\in\left\{  1,2,\ldots,M\right\}  $, the element $%
\begin{cases}
s_{i}, & \text{if }i<M;\\
t, & \text{if }i=M
\end{cases}
$ is a well-defined element of $Z_{i}$\ \ \ \ \footnote{\textit{Proof.} Let
$i\in\left\{  1,2,\ldots,M\right\}  $. We must prove that
\begin{equation}%
\begin{cases}
s_{i}, & \text{if }i<M;\\
t, & \text{if }i=M
\end{cases}
\text{ is a well-defined element of }Z_{i}.
\label{pf.lem.prodrule.S.Psi.pf.fn1.1}%
\end{equation}
\par
We are in one of the following two cases:
\par
\textit{Case 1:} We have $i\neq M$.
\par
\textit{Case 2:} We have $i=M$.
\par
Let us first consider Case 1. In this case, we have $i\neq M$. Combining this
with $i\in\left\{  1,2,\ldots,M\right\}  $, we obtain $i\in\left\{
1,2,\ldots,M\right\}  \setminus\left\{  M\right\}  =\left\{  1,2,\ldots
,M-1\right\}  $. Hence, $s_{i}\in Z_{i}$ (by (\ref{pf.lem.prodrule.S.Psi.pf.1}%
)). Thus, $s_{i}$ is a well-defined element of $Z_{i}$. In other words, $%
\begin{cases}
s_{i}, & \text{if }i<M;\\
t, & \text{if }i=M
\end{cases}
$ is a well-defined element of $Z_{i}$ (since $%
\begin{cases}
s_{i}, & \text{if }i<M;\\
t, & \text{if }i=M
\end{cases}
=s_{i}$ (since $i<M$ (since $i\in\left\{  1,2,\ldots,M-1\right\}  $))). Thus,
(\ref{pf.lem.prodrule.S.Psi.pf.fn1.1}) is proven in Case 1.
\par
Let us now consider Case 2. In this case, we have $i=M$. But $t$ is a
well-defined element of $Z_{M}$. In other words, $%
\begin{cases}
s_{i}, & \text{if }i<M;\\
t, & \text{if }i=M
\end{cases}
$ is a well-defined element of $Z_{i}$ (since $%
\begin{cases}
s_{i}, & \text{if }i<M;\\
t, & \text{if }i=M
\end{cases}
=t$ (since $i=M$) and $i=M$). Thus, (\ref{pf.lem.prodrule.S.Psi.pf.fn1.1}) is
proven in Case 2.
\par
We have now proven (\ref{pf.lem.prodrule.S.Psi.pf.fn1.1}) in each of the two
Cases 1 and 2. Thus, (\ref{pf.lem.prodrule.S.Psi.pf.fn1.1}) always holds
(since Cases 1 and 2 cover all possibilities). In other words, $%
\begin{cases}
s_{i}, & \text{if }i<M;\\
t, & \text{if }i=M
\end{cases}
$ is a well-defined element of $Z_{i}$. Qed.}. Hence, we can define an
$M$-tuple $\left(  q_{1},q_{2},\ldots,q_{M}\right)  \in Z_{1}\times
Z_{2}\times\cdots\times Z_{M}$ by%
\[
\left(  q_{i}=%
\begin{cases}
s_{i}, & \text{if }i<M;\\
t, & \text{if }i=M
\end{cases}
\ \ \ \ \ \ \ \ \ \ \text{for every }i\in\left\{  1,2,\ldots,M\right\}
\right)  .
\]
Consider this $M$-tuple $\left(  q_{1},q_{2},\ldots,q_{M}\right)  $. The
definition of this $M$-tuple shows that%
\[
\left(  q_{1},q_{2},\ldots,q_{M}\right)  =\left(  s_{1},s_{2},\ldots
,s_{M-1},t\right)  .
\]

For every $i\in\left\{  1,2,\ldots,M-1\right\}  $, we have%
\[
q_{i}=%
\begin{cases}
s_{i}, & \text{if }i<M;\\
t, & \text{if }i=M
\end{cases}
=s_{i}\ \ \ \ \ \ \ \ \ \ \left(  \text{since }i<M\text{ (since }i\in\left\{
1,2,\ldots,M-1\right\}  \text{)}\right)  .
\]
In other words, $\left(  q_{1},q_{2},\ldots,q_{M-1}\right)  =\left(
s_{1},s_{2},\ldots,s_{M-1}\right)  $. Also, the definition of $q_{M}$ shows
that $q_{M}=%
\begin{cases}
s_{M}, & \text{if }M<M;\\
t, & \text{if }M=M
\end{cases}
=t$ (since $M=M$).

But the definition of $\Psi$ yields
\[
\Psi\left(  q_{1},q_{2},\ldots,q_{M}\right)  =\left(  \underbrace{\left(
q_{1},q_{2},\ldots,q_{M-1}\right)  }_{=\left(  s_{1},s_{2},\ldots
,s_{M-1}\right)  },\underbrace{q_{M}}_{=t}\right)  =\left(  \left(
s_{1},s_{2},\ldots,s_{M-1}\right)  ,t\right)  .
\]
Comparing this with $\Psi\underbrace{\left(  q_{1},q_{2},\ldots,q_{M}\right)
}_{=\left(  s_{1},s_{2},\ldots,s_{M-1},t\right)  }=\Psi\left(  s_{1}%
,s_{2},\ldots,s_{M-1},t\right)  $, we obtain
\[
\Psi\left(  s_{1},s_{2},\ldots,s_{M-1},t\right)  =\left(  \left(  s_{1}%
,s_{2},\ldots,s_{M-1}\right)  ,t\right)  .
\]
This proves (\ref{pf.lem.prodrule.S.Psi}).]

Now, we have $\Phi\circ\Psi=\operatorname*{id}$%
\ \ \ \ \footnote{\textit{Proof.} Let $\alpha\in Z_{1}\times Z_{2}\times
\cdots\times Z_{M}$. Thus, we can write $\alpha$ in the form $\left(
s_{1},s_{2},\ldots,s_{M}\right)  $ for some $\left(  s_{1},s_{2},\ldots
,s_{M}\right)  \in Z_{1}\times Z_{2}\times\cdots\times Z_{M}$. Consider this
$\left(  s_{1},s_{2},\ldots,s_{M}\right)  $. Hence, $\alpha=\left(
s_{1},s_{2},\ldots,s_{M}\right)  $.
\par
Now,%
\begin{align*}
\left(  \Phi\circ\Psi\right)  \left(  \underbrace{\alpha}_{=\left(
s_{1},s_{2},\ldots,s_{M}\right)  }\right)   &  =\left(  \Phi\circ\Psi\right)
\left(  s_{1},s_{2},\ldots,s_{M}\right)  =\Phi\left(  \underbrace{\Psi\left(
s_{1},s_{2},\ldots,s_{M}\right)  }_{\substack{=\left(  \left(  s_{1}%
,s_{2},\ldots,s_{M-1}\right)  ,s_{M}\right)  \\\text{(by the definition of
}\Psi\text{)}}}\right) \\
&  =\Phi\left(  \left(  s_{1},s_{2},\ldots,s_{M-1}\right)  ,s_{M}\right)
=\left(  s_{1},s_{2},\ldots,s_{M-1},s_{M}\right) \\
&  \ \ \ \ \ \ \ \ \ \ \left(  \text{by the definition of }\Phi\right) \\
&  =\left(  s_{1},s_{2},\ldots,s_{M}\right)  =\alpha=\operatorname*{id}\left(
\alpha\right)  .
\end{align*}
\par
Let us now forget that we fixed $\alpha$. We thus have shown that $\left(
\Phi\circ\Psi\right)  \left(  \alpha\right)  =\operatorname*{id}\left(
\alpha\right)  $ for every $\alpha\in Z_{1}\times Z_{2}\times\cdots\times
Z_{M}$. In other words, $\Phi\circ\Psi=\operatorname*{id}$, qed.} and
$\Psi\circ\Phi=\operatorname*{id}$\ \ \ \ \footnote{\textit{Proof.} Let
$\beta\in\left(  Z_{1}\times Z_{2}\times\cdots\times Z_{M-1}\right)  \times
Z_{M}$. Thus, we can write $\beta$ in the form $\left(  \gamma,t\right)  $ for
some $\gamma\in Z_{1}\times Z_{2}\times\cdots\times Z_{M-1}$ and $t\in Z_{M}$.
Consider these $\gamma$ and $t$. Thus, $\beta=\left(  \gamma,t\right)  $.
\par
We have $\gamma\in Z_{1}\times Z_{2}\times\cdots\times Z_{M-1}$. Thus, we can
write $\gamma$ in the form $\left(  s_{1},s_{2},\ldots,s_{M-1}\right)  $ for
some $\left(  s_{1},s_{2},\ldots,s_{M-1}\right)  \in Z_{1}\times Z_{2}%
\times\cdots\times Z_{M-1}$. Consider this $\left(  s_{1},s_{2},\ldots
,s_{M-1}\right)  $. Hence, $\gamma=\left(  s_{1},s_{2},\ldots,s_{M-1}\right)
$.
\par
Now, $\beta=\left(  \underbrace{\gamma}_{=\left(  s_{1},s_{2},\ldots
,s_{M-1}\right)  },t\right)  =\left(  \left(  s_{1},s_{2},\ldots
,s_{M-1}\right)  ,t\right)  $. Applying the map $\Psi\circ\Phi$ to both sides
of this equality, we obtain%
\begin{align*}
\left(  \Psi\circ\Phi\right)  \left(  \beta\right)   &  =\left(  \Psi\circ
\Phi\right)  \left(  \left(  s_{1},s_{2},\ldots,s_{M-1}\right)  ,t\right)
=\Psi\left(  \underbrace{\Phi\left(  \left(  s_{1},s_{2},\ldots,s_{M-1}%
\right)  ,t\right)  }_{\substack{=\left(  s_{1},s_{2},\ldots,s_{M-1},t\right)
\\\text{(by the definition of }\Phi\text{)}}}\right) \\
&  =\Psi\left(  s_{1},s_{2},\ldots,s_{M-1},t\right)  =\left(  \left(
s_{1},s_{2},\ldots,s_{M-1}\right)  ,t\right)  \ \ \ \ \ \ \ \ \ \ \left(
\text{by (\ref{pf.lem.prodrule.S.Psi})}\right) \\
&  =\beta=\operatorname*{id}\left(  \beta\right)  .
\end{align*}
\par
Let us now forget that we fixed $\beta$. We thus have shown that $\left(
\Psi\circ\Phi\right)  \left(  \beta\right)  =\operatorname*{id}\left(
\beta\right)  $ for every $\beta\in\left(  Z_{1}\times Z_{2}\times\cdots\times
Z_{M-1}\right)  \times Z_{M}$. In other words, $\Psi\circ\Phi
=\operatorname*{id}$, qed.}. These two equalities show that the maps $\Phi$
and $\Psi$ are mutually inverse. Thus, the map $\Psi$ is invertible. In other
words, the map $\Psi$ is a bijection.

Now, $\Psi$ is the map
\begin{align*}
Z_{1}\times Z_{2}\times\cdots\times Z_{M}  &  \rightarrow\left(  Z_{1}\times
Z_{2}\times\cdots\times Z_{M-1}\right)  \times Z_{M},\\
\left(  s_{1},s_{2},\ldots,s_{M}\right)   &  \mapsto\left(  \left(
s_{1},s_{2},\ldots,s_{M-1}\right)  ,s_{M}\right)
\end{align*}
(because $\Psi$ is the map $Z_{1}\times Z_{2}\times\cdots\times Z_{M}%
\rightarrow\left(  Z_{1}\times Z_{2}\times\cdots\times Z_{M-1}\right)  \times
Z_{M}$ satisfying (\ref{pf.lem.prodrule.S.Psidef})). Thus, the map
\begin{align*}
Z_{1}\times Z_{2}\times\cdots\times Z_{M}  &  \rightarrow\left(  Z_{1}\times
Z_{2}\times\cdots\times Z_{M-1}\right)  \times Z_{M},\\
\left(  s_{1},s_{2},\ldots,s_{M}\right)   &  \mapsto\left(  \left(
s_{1},s_{2},\ldots,s_{M-1}\right)  ,s_{M}\right)
\end{align*}
is a bijection (since $\Psi$ is a bijection). This proves Lemma
\ref{lem.prodrule.prod-assM}.
\end{proof}
\end{verlong}

For the sake of convenience, let us state a particular case of Lemma
\ref{lem.prodrule.prod-assM}:

\begin{corollary}
\label{cor.prodrule.prod-assMZ}Let $M$ be a positive integer. Let $Z$ be a
set. Then, the map
\begin{align*}
Z^{M}  &  \rightarrow Z^{M-1}\times Z,\\
\left(  s_{1},s_{2},\ldots,s_{M}\right)   &  \mapsto\left(  \left(
s_{1},s_{2},\ldots,s_{M-1}\right)  ,s_{M}\right)
\end{align*}
is a bijection.
\end{corollary}

\begin{proof}
[Proof of Corollary \ref{cor.prodrule.prod-assMZ}.]Lemma
\ref{lem.prodrule.prod-assM} (applied to $Z_{i}=Z$) shows that the map%
\begin{align*}
\underbrace{Z\times Z\times\cdots\times Z}_{M\text{ factors}}  &
\rightarrow\left(  \underbrace{Z\times Z\times\cdots\times Z}_{M-1\text{
factors}}\right)  \times Z,\\
\left(  s_{1},s_{2},\ldots,s_{M}\right)   &  \mapsto\left(  \left(
s_{1},s_{2},\ldots,s_{M-1}\right)  ,s_{M}\right)
\end{align*}
is a bijection. Since $\underbrace{Z\times Z\times\cdots\times Z}_{M\text{
factors}}=Z^{M}$ and $\underbrace{Z\times Z\times\cdots\times Z}_{M-1\text{
factors}}=Z^{M-1}$, this can be rewritten as follows: The map%
\begin{align*}
Z^{M}  &  \rightarrow Z^{M-1}\times Z,\\
\left(  s_{1},s_{2},\ldots,s_{M}\right)   &  \mapsto\left(  \left(
s_{1},s_{2},\ldots,s_{M-1}\right)  ,s_{M}\right)
\end{align*}
is a bijection. This proves Corollary \ref{cor.prodrule.prod-assMZ}.
\end{proof}

Next, we state a lemma that is essentially the statement of Exercise
\ref{exe.multichoose}:

\begin{lemma}
\label{lem.sol.multichoose.sum}Every $n\in\mathbb{N}$ and $m\in\mathbb{N}$
satisfy%
\[
\sum_{\substack{\left(  a_{1},a_{2},\ldots,a_{m}\right)  \in\mathbb{N}%
^{m};\\a_{1}+a_{2}+\cdots+a_{m}=n}}1=\dbinom{n+m-1}{n}.
\]

\end{lemma}

\begin{proof}
[Proof of Lemma \ref{lem.sol.multichoose.sum}.]We shall prove Lemma
\ref{lem.sol.multichoose.sum} by induction over $m$:

\begin{vershort}
\textit{Induction base:} Lemma \ref{lem.sol.multichoose.sum} holds for
$m=0$\ \ \ \ \footnote{\textit{Proof.} We must show that Lemma
\ref{lem.sol.multichoose.sum} holds for $m=0$. In other words, we must prove
that $\sum_{\substack{\left(  a_{1},a_{2},\ldots,a_{0}\right)  \in
\mathbb{N}^{0};\\a_{1}+a_{2}+\cdots+a_{0}=n}}1=\dbinom{n+0-1}{n}$ for all
$n\in\mathbb{N}$.
\par
Fix $n\in\mathbb{N}$. We are in one of the following two cases:
\par
\textit{Case 1:} We have $n=0$.
\par
\textit{Case 2:} We have $n\neq0$.
\par
Let us first consider Case 1. In this case, we have $n=0$. There exists
exactly one $0$-tuple $\left(  a_{1},a_{2},\ldots,a_{0}\right)  \in
\mathbb{N}^{0}$ (namely, the empty list $\left(  {}\right)  $), and this
$0$-tuple satisfies $a_{1}+a_{2}+\cdots+a_{0}=n$ (because it satisfies
$a_{1}+a_{2}+\cdots+a_{0}=\left(  \text{empty sum}\right)  =0=n$). Thus, the
sum $\sum_{\substack{\left(  a_{1},a_{2},\ldots,a_{0}\right)  \in
\mathbb{N}^{0};\\a_{1}+a_{2}+\cdots+a_{0}=n}}1$ has exactly one addend
(namely, the addend corresponding to $\left(  a_{1},a_{2},\ldots,a_{0}\right)
=\left(  {}\right)  $). Hence, this sum rewrites as follows:%
\[
\sum_{\substack{\left(  a_{1},a_{2},\ldots,a_{0}\right)  \in\mathbb{N}%
^{0};\\a_{1}+a_{2}+\cdots+a_{0}=n}}1=1.
\]
Comparing this with%
\begin{align*}
\dbinom{n+0-1}{n}  &  =\dbinom{0+0-1}{0}\ \ \ \ \ \ \ \ \ \ \left(
\text{since }n=0\right) \\
&  =1\ \ \ \ \ \ \ \ \ \ \left(  \text{by Proposition \ref{prop.binom.00}
\textbf{(a)} (applied to }0+0-1\text{ instead of }m\text{)}\right)  ,
\end{align*}
we obtain $\sum_{\substack{\left(  a_{1},a_{2},\ldots,a_{0}\right)
\in\mathbb{N}^{0};\\a_{1}+a_{2}+\cdots+a_{0}=n}}1=\dbinom{n+0-1}{n}$. Hence,
the equality $\sum_{\substack{\left(  a_{1},a_{2},\ldots,a_{0}\right)
\in\mathbb{N}^{0};\\a_{1}+a_{2}+\cdots+a_{0}=n}}1=\dbinom{n+0-1}{n}$ is proven
in Case 1.
\par
Let us now consider Case 2. In this case, we have $n\neq0$. Hence, $n$ is a
positive integer (since $n\in\mathbb{N}$). Thus, $n-1\in\mathbb{N}$.
\par
Each $0$-tuple $\left(  a_{1},a_{2},\ldots,a_{0}\right)  \in\mathbb{N}^{0}$
satisfies $a_{1}+a_{2}+\cdots+a_{0}=\left(  \text{empty sum}\right)  =0\neq
n$. In other words, no $\left(  a_{1},a_{2},\ldots,a_{0}\right)  \in
\mathbb{N}^{0}$ satisfies $a_{1}+a_{2}+\cdots+a_{0}=n$. Hence, the sum
$\sum_{\substack{\left(  a_{1},a_{2},\ldots,a_{0}\right)  \in\mathbb{N}%
^{0};\\a_{1}+a_{2}+\cdots+a_{0}=n}}1$ is an empty sum. Therefore,%
\[
\sum_{\substack{\left(  a_{1},a_{2},\ldots,a_{0}\right)  \in\mathbb{N}%
^{0};\\a_{1}+a_{2}+\cdots+a_{0}=n}}1=\left(  \text{empty sum}\right)  =0.
\]
\par
But $n-1<n$ and $n-1\in\mathbb{N}$. Thus, Proposition \ref{prop.binom.0}
(applied to $n-1$ instead of $m$) yields $\dbinom{n-1}{n}=0$. Hence,
$\dbinom{n+0-1}{n}=\dbinom{n-1}{n}=0$. Thus, $\sum_{\substack{\left(
a_{1},a_{2},\ldots,a_{0}\right)  \in\mathbb{N}^{0};\\a_{1}+a_{2}+\cdots
+a_{0}=n}}1=0=\dbinom{n+0-1}{n}$. Hence, the equality $\sum_{\substack{\left(
a_{1},a_{2},\ldots,a_{0}\right)  \in\mathbb{N}^{0};\\a_{1}+a_{2}+\cdots
+a_{0}=n}}1=\dbinom{n+0-1}{n}$ is proven in Case 2.
\par
We have now proven this equality in each of the two Cases 1 and 2. Thus, this
equality always holds. In other words, Lemma \ref{lem.sol.multichoose.sum}
holds for $m=0$.}. This completes the induction base.
\end{vershort}

\begin{verlong}
\textit{Induction base:} Lemma \ref{lem.sol.multichoose.sum} holds for
$m=0$\ \ \ \ \footnote{\textit{Proof.} Assume that $m=0$. We must then show
that Lemma \ref{lem.sol.multichoose.sum} holds.
\par
We are in one of the following two cases:
\par
\textit{Case 1:} We have $n=0$.
\par
\textit{Case 2:} We have $n\neq0$.
\par
Let us first consider Case 1. In this case, we have $n=0$. There exists
exactly one $0$-tuple $\left(  a_{1},a_{2},\ldots,a_{0}\right)  \in
\mathbb{N}^{0}$ (namely, the empty list $\left(  {}\right)  $). Hence, the sum
$\sum_{\left(  a_{1},a_{2},\ldots,a_{0}\right)  \in\mathbb{N}^{0}}1$ has
exactly one addend (namely, the addend corresponding to $\left(  a_{1}%
,a_{2},\ldots,a_{0}\right)  =\left(  {}\right)  $). Hence, this sum rewrites
as follows:%
\[
\sum_{\left(  a_{1},a_{2},\ldots,a_{0}\right)  \in\mathbb{N}^{0}}1=1.
\]
But recall that $m=0$. Hence,%
\[
\sum_{\substack{\left(  a_{1},a_{2},\ldots,a_{m}\right)  \in\mathbb{N}%
^{m};\\a_{1}+a_{2}+\cdots+a_{m}=n}}1=\underbrace{\sum_{\substack{\left(
a_{1},a_{2},\ldots,a_{0}\right)  \in\mathbb{N}^{0};\\a_{1}+a_{2}+\cdots
+a_{0}=n}}}_{\substack{=\sum_{\left(  a_{1},a_{2},\ldots,a_{0}\right)
\in\mathbb{N}^{0}}\\\text{(since each }\left(  a_{1},a_{2},\ldots
,a_{0}\right)  \in\mathbb{N}^{0}\text{ automatically}\\\text{satisfies }%
a_{1}+a_{2}+\cdots+a_{0}=n\\\text{(because for each }\left(  a_{1}%
,a_{2},\ldots,a_{0}\right)  \in\mathbb{N}^{0}\text{,}\\\text{we have }%
a_{1}+a_{2}+\cdots+a_{0}=\left(  \text{empty sum}\right)  =0=n\text{))}%
}}1=\sum_{\left(  a_{1},a_{2},\ldots,a_{0}\right)  \in\mathbb{N}^{0}}1=1.
\]
Comparing this with%
\begin{align*}
\dbinom{n+m-1}{n}  &  =\dbinom{0+m-1}{0}\ \ \ \ \ \ \ \ \ \ \left(
\text{since }n=0\right) \\
&  =1\ \ \ \ \ \ \ \ \ \ \left(  \text{by Proposition \ref{prop.binom.00}
\textbf{(a)} (applied to }0+m-1\text{ instead of }m\text{)}\right)  ,
\end{align*}
we obtain%
\[
\sum_{\substack{\left(  a_{1},a_{2},\ldots,a_{m}\right)  \in\mathbb{N}%
^{m};\\a_{1}+a_{2}+\cdots+a_{m}=n}}1=\dbinom{n+m-1}{n}.
\]
Hence, Lemma \ref{lem.sol.multichoose.sum} holds. We thus have proved Lemma
\ref{lem.sol.multichoose.sum} in Case 1.
\par
Let us now consider Case 2. In this case, we have $n\neq0$. Hence, $n$ is a
positive integer (since $n\in\mathbb{N}$). Thus, $n-1\in\mathbb{N}$.
\par
Each $0$-tuple $\left(  a_{1},a_{2},\ldots,a_{0}\right)  \in\mathbb{N}^{0}$
satisfies $a_{1}+a_{2}+\cdots+a_{0}=\left(  \text{empty sum}\right)  =0\neq
n$. In other words, no $\left(  a_{1},a_{2},\ldots,a_{0}\right)  \in
\mathbb{N}^{0}$ satisfies $a_{1}+a_{2}+\cdots+a_{0}=n$. Hence, the sum
$\sum_{\substack{\left(  a_{1},a_{2},\ldots,a_{0}\right)  \in\mathbb{N}%
^{0};\\a_{1}+a_{2}+\cdots+a_{0}=n}}1$ is an empty sum. Therefore, this sum
rewrites as follows:%
\[
\sum_{\substack{\left(  a_{1},a_{2},\ldots,a_{0}\right)  \in\mathbb{N}%
^{0};\\a_{1}+a_{2}+\cdots+a_{0}=n}}1=\left(  \text{empty sum}\right)  =0.
\]
\par
Now, recall that $m=0$. Also, $n-1<n$ and $n-1\in\mathbb{N}$. Thus,
Proposition \ref{prop.binom.0} (applied to $n-1$ instead of $m$) yields
$\dbinom{n-1}{n}=0$. But from $m=0$, we obtain $n+m-1=n+0-1=n-1$; therefore,
$\dbinom{n+m-1}{n}=\dbinom{n-1}{n}=0$.
\par
From $m=0$, we obtain%
\[
\sum_{\substack{\left(  a_{1},a_{2},\ldots,a_{m}\right)  \in\mathbb{N}%
^{m};\\a_{1}+a_{2}+\cdots+a_{m}=n}}1=\sum_{\substack{\left(  a_{1}%
,a_{2},\ldots,a_{0}\right)  \in\mathbb{N}^{0};\\a_{1}+a_{2}+\cdots+a_{0}%
=n}}1=0=\dbinom{n+m-1}{n}%
\]
(since $\dbinom{n+m-1}{n}=0$). Hence, Lemma \ref{lem.sol.multichoose.sum}
holds. We thus have proved Lemma \ref{lem.sol.multichoose.sum} in Case 2.
\par
We now have proved Lemma \ref{lem.sol.multichoose.sum} in each of the two
Cases 1 and 2. Thus, Lemma \ref{lem.sol.multichoose.sum} always holds (under
the assumption that $m=0$). Qed.}. This completes the induction base.
\end{verlong}

\textit{Induction step:} Fix a positive integer $M$. Assume that Lemma
\ref{lem.sol.multichoose.sum} holds for $m=M-1$. We now must show that Lemma
\ref{lem.sol.multichoose.sum} holds for $m=M$.

We have assumed that Lemma \ref{lem.sol.multichoose.sum} holds for $m=M-1$. In
other words, every $n\in\mathbb{N}$ satisfies%
\begin{equation}
\sum_{\substack{\left(  a_{1},a_{2},\ldots,a_{M-1}\right)  \in\mathbb{N}%
^{M-1};\\a_{1}+a_{2}+\cdots+a_{M-1}=n}}1=\dbinom{n+\left(  M-1\right)  -1}{n}.
\label{pf.lem.sol.multichoose.sum.IH}%
\end{equation}

We must show that Lemma \ref{lem.sol.multichoose.sum} holds for $m=M$. In
other words, we must prove that every $n\in\mathbb{N}$ satisfies%
\begin{equation}
\sum_{\substack{\left(  a_{1},a_{2},\ldots,a_{M}\right)  \in\mathbb{N}%
^{M};\\a_{1}+a_{2}+\cdots+a_{M}=n}}1=\dbinom{n+M-1}{n}.
\label{pf.lem.sol.multichoose.sum.IG}%
\end{equation}
Let us now prove this.

Let $n\in\mathbb{N}$. For each $\left(  a_{1},a_{2},\ldots,a_{M}\right)
\in\mathbb{N}^{M}$ satisfying $a_{1}+a_{2}+\cdots+a_{M}=n$, we have
$n-a_{M}\in\left\{  0,1,\ldots,n\right\}  $\ \ \ \ \footnote{\textit{Proof.}
Let $\left(  a_{1},a_{2},\ldots,a_{M}\right)  \in\mathbb{N}^{M}$ be such that
$a_{1}+a_{2}+\cdots+a_{M}=n$. We must show that $n-a_{M}\in\left\{
0,1,\ldots,n\right\}  $.
\par
We have $\left(  a_{1},a_{2},\ldots,a_{M}\right)  \in\mathbb{N}^{M}$. Thus,
$a_{1},a_{2},\ldots,a_{M}$ are elements of $\mathbb{N}$. Hence, in particular,
$a_{1},a_{2},\ldots,a_{M-1}$ are elements of $\mathbb{N}$. Thus, $a_{1}%
+a_{2}+\cdots+a_{M-1}\in\mathbb{N}$, so that $a_{1}+a_{2}+\cdots+a_{M-1}\geq
0$.
\par
From $a_{1}+a_{2}+\cdots+a_{M}=n$, we obtain%
\begin{align*}
n  &  =a_{1}+a_{2}+\cdots+a_{M}=\underbrace{\left(  a_{1}+a_{2}+\cdots
+a_{M-1}\right)  }_{\geq0}+a_{M}\ \ \ \ \ \ \ \ \ \ \left(  \text{since
}M\text{ is a positive integer}\right) \\
&  \geq a_{M},
\end{align*}
so that $n-a_{M}\geq0$. Also, we know that $a_{1},a_{2},\ldots,a_{M}$ are
elements of $\mathbb{N}$; thus, $a_{M}\in\mathbb{N}$. Hence, $a_{M}\geq0$.
\par
Clearly, $n$ and $a_{M}$ are integers (since $n\in\mathbb{N}\subseteq
\mathbb{Z}$ and $a_{M}\in\mathbb{N}\subseteq\mathbb{Z}$). Thus, $n-a_{M}$ is
an integer. Combining this with $n-a_{M}\geq0$, we conclude that $n-a_{M}%
\in\mathbb{N}$. Combining this with $n-\underbrace{a_{M}}_{\geq0}\leq n$, we
obtain $n-a_{M}\in\left\{  0,1,\ldots,n\right\}  $. Qed.}. Hence, we have the
following equality of summation signs:%
\begin{equation}
\sum_{\substack{\left(  a_{1},a_{2},\ldots,a_{M}\right)  \in\mathbb{N}%
^{M};\\a_{1}+a_{2}+\cdots+a_{M}=n}}=\underbrace{\sum_{r\in\left\{
0,1,\ldots,n\right\}  }}_{=\sum_{r=0}^{n}}\sum_{\substack{\left(  a_{1}%
,a_{2},\ldots,a_{M}\right)  \in\mathbb{N}^{M};\\a_{1}+a_{2}+\cdots
+a_{M}=n;\\n-a_{M}=r}}=\sum_{r=0}^{n}\sum_{\substack{\left(  a_{1}%
,a_{2},\ldots,a_{M}\right)  \in\mathbb{N}^{M};\\a_{1}+a_{2}+\cdots
+a_{M}=n;\\n-a_{M}=r}}. \label{pf.lem.sol.multichoose.sum.sums1}%
\end{equation}

Now, we are going to show that%
\begin{equation}
\sum_{\substack{\left(  a_{1},a_{2},\ldots,a_{M}\right)  \in\mathbb{N}%
^{M};\\a_{1}+a_{2}+\cdots+a_{M}=n;\\n-a_{M}=r}}1=\dbinom{r+M-2}{r}
\label{pf.lem.sol.multichoose.sum.sums2}%
\end{equation}
for each $r\in\left\{  0,1,\ldots,n\right\}  $.

[\textit{Proof of (\ref{pf.lem.sol.multichoose.sum.sums2}):} Let $r\in\left\{
0,1,\ldots,n\right\}  $.

\begin{vershort}
For any $\left(  a_{1},a_{2},\ldots,a_{M}\right)  \in\mathbb{N}^{M}$, the
condition $\left(  a_{1}+a_{2}+\cdots+a_{M}=n\text{ and }n-a_{M}=r\right)  $
is equivalent to the condition $\left(  a_{1}+a_{2}+\cdots+a_{M-1}=r\text{ and
}a_{M}=n-r\right)  $.\ \ \ \ \footnote{This can be checked easily. For
example, in order to prove the implication%
\[
\left(  a_{1}+a_{2}+\cdots+a_{M}=n\text{ and }n-a_{M}=r\right)
\ \Longrightarrow\ \left(  a_{1}+a_{2}+\cdots+a_{M-1}=r\text{ and }%
a_{M}=n-r\right)  ,
\]
it suffices to assume that $\left(  a_{1}+a_{2}+\cdots+a_{M}=n\text{ and
}n-a_{M}=r\right)  $ holds, and then to conclude that%
\begin{align*}
a_{1}+a_{2}+\cdots+a_{M-1}  &  =\underbrace{\left(  a_{1}+a_{2}+\cdots
+a_{M}\right)  }_{=n}-a_{M}=n-a_{M}=r\ \ \ \ \ \ \ \ \ \ \text{and}\\
a_{M}  &  =n-r\ \ \ \ \ \ \ \ \ \ \left(  \text{since }n-a_{M}=r\right)  .
\end{align*}
The converse implication is proven similarly.}
\end{vershort}

\begin{verlong}
Fix $\left(  a_{1},a_{2},\ldots,a_{M}\right)  \in\mathbb{N}^{M}$. We shall
show that the condition \newline$\left(  a_{1}+a_{2}+\cdots+a_{M}=n\text{ and
}n-a_{M}=r\right)  $ is equivalent to the condition \newline$\left(
a_{1}+a_{2}+\cdots+a_{M-1}=r\text{ and }a_{M}=n-r\right)  $.

Indeed, we have the logical implication%
\begin{align}
&  \ \ \left(  a_{1}+a_{2}+\cdots+a_{M}=n\text{ and }n-a_{M}=r\right)
\nonumber\\
&  \Longrightarrow\ \ \left(  a_{1}+a_{2}+\cdots+a_{M-1}=r\text{ and }%
a_{M}=n-r\right)  \label{pf.lem.sol.multichoose.sum.sums2.pf.i1}%
\end{align}
\footnote{\textit{Proof of (\ref{pf.lem.sol.multichoose.sum.sums2.pf.i1}):}
Assume that $\left(  a_{1}+a_{2}+\cdots+a_{M}=n\text{ and }n-a_{M}=r\right)
$. Then, $a_{1}+a_{2}+\cdots+a_{M}=n$, so that%
\[
n=a_{1}+a_{2}+\cdots+a_{M}=\left(  a_{1}+a_{2}+\cdots+a_{M-1}\right)
+a_{M}\ \ \ \ \ \ \ \ \ \ \left(  \text{since }M\text{ is a positive
integer}\right)  .
\]
Subtracting $a_{M}$ from both sides of this equality, we obtain $n-a_{M}%
=a_{1}+a_{2}+\cdots+a_{M-1}$. Comparing this with $n-a_{M}=r$, we find
$a_{1}+a_{2}+\cdots+a_{M-1}=r$. Furthermore, from $n-a_{M}=r$, we obtain
$a_{M}=n-r$. We thus have proven that $\left(  a_{1}+a_{2}+\cdots
+a_{M-1}=r\text{ and }a_{M}=n-r\right)  $.
\par
Now, forget our assumption that $\left(  a_{1}+a_{2}+\cdots+a_{M}=n\text{ and
}n-a_{M}=r\right)  $. We thus have shown that $\left(  a_{1}+a_{2}%
+\cdots+a_{M-1}=r\text{ and }a_{M}=n-r\right)  $ holds under the assumption
that $\left(  a_{1}+a_{2}+\cdots+a_{M}=n\text{ and }n-a_{M}=r\right)  $. In
other words, we have proven the implication
(\ref{pf.lem.sol.multichoose.sum.sums2.pf.i1}).} and the logical implication%
\begin{align}
&  \ \ \left(  a_{1}+a_{2}+\cdots+a_{M-1}=r\text{ and }a_{M}=n-r\right)
\nonumber\\
&  \Longrightarrow\ \ \left(  a_{1}+a_{2}+\cdots+a_{M}=n\text{ and }%
n-a_{M}=r\right)  \label{pf.lem.sol.multichoose.sum.sums2.pf.i2}%
\end{align}
\footnote{\textit{Proof of (\ref{pf.lem.sol.multichoose.sum.sums2.pf.i2}):}
Assume that $\left(  a_{1}+a_{2}+\cdots+a_{M-1}=r\text{ and }a_{M}=n-r\right)
$. Then, $n-a_{M}=r$ (since $a_{M}=n-r$) and
\begin{align*}
a_{1}+a_{2}+\cdots+a_{M}  &  =\underbrace{\left(  a_{1}+a_{2}+\cdots
+a_{M-1}\right)  }_{=r}+a_{M}\ \ \ \ \ \ \ \ \ \ \left(  \text{since }M\text{
is a positive integer}\right) \\
&  =r+a_{M}=n\ \ \ \ \ \ \ \ \ \ \left(  \text{since }a_{M}=n-r\right)  .
\end{align*}
We thus have proven that $\left(  a_{1}+a_{2}+\cdots+a_{M}=n\text{ and
}n-a_{M}=r\right)  $.
\par
Now, forget our assumption that $\left(  a_{1}+a_{2}+\cdots+a_{M-1}=r\text{
and }a_{M}=n-r\right)  $. We thus have shown that $\left(  a_{1}+a_{2}%
+\cdots+a_{M}=n\text{ and }n-a_{M}=r\right)  $ holds under the assumption that
$\left(  a_{1}+a_{2}+\cdots+a_{M-1}=r\text{ and }a_{M}=n-r\right)  $. In other
words, we have proven the implication
(\ref{pf.lem.sol.multichoose.sum.sums2.pf.i2}).}. Combining these two
implications, we obtain the logical equivalence%
\begin{align*}
&  \ \ \left(  a_{1}+a_{2}+\cdots+a_{M}=n\text{ and }n-a_{M}=r\right) \\
&  \Longleftrightarrow\ \ \left(  a_{1}+a_{2}+\cdots+a_{M-1}=r\text{ and
}a_{M}=n-r\right)  .
\end{align*}
In other words, the condition $\left(  a_{1}+a_{2}+\cdots+a_{M}=n\text{ and
}n-a_{M}=r\right)  $ is equivalent to the condition $\left(  a_{1}%
+a_{2}+\cdots+a_{M-1}=r\text{ and }a_{M}=n-r\right)  $.

Now, forget that we fixed $\left(  a_{1},a_{2},\ldots,a_{M}\right)  $. We thus
have shown that for any $\left(  a_{1},a_{2},\ldots,a_{M}\right)
\in\mathbb{N}^{M}$, the condition $\left(  a_{1}+a_{2}+\cdots+a_{M}=n\text{
and }n-a_{M}=r\right)  $ is equivalent to the condition $\left(  a_{1}%
+a_{2}+\cdots+a_{M-1}=r\text{ and }a_{M}=n-r\right)  $.
\end{verlong}

Thus, we have the following equality of summation signs:%
\begin{equation}
\sum_{\substack{\left(  a_{1},a_{2},\ldots,a_{M}\right)  \in\mathbb{N}%
^{M};\\a_{1}+a_{2}+\cdots+a_{M}=n;\\n-a_{M}=r}}=\sum_{\substack{\left(
a_{1},a_{2},\ldots,a_{M}\right)  \in\mathbb{N}^{M};\\a_{1}+a_{2}%
+\cdots+a_{M-1}=r;\\a_{M}=n-r}}. \label{pf.lem.sol.multichoose.sum.sums3}%
\end{equation}

We have $r\leq n$ (since $r\in\left\{  0,1,\ldots,n\right\}  $), thus $n\geq
r$ and therefore $n-r\geq0$. Hence, $n-r\in\mathbb{N}$. Thus, the sum
$\sum_{\substack{q\in\mathbb{N};\\q=n-r}}1$ has exactly one addend (namely,
the addend for $q=n-r$). Hence, it simplifies as follows:%
\begin{equation}
\sum_{\substack{q\in\mathbb{N};\\q=n-r}}1=1.
\label{pf.lem.sol.multichoose.sum.sums4}%
\end{equation}

Corollary \ref{cor.prodrule.prod-assMZ} (applied to $Z=\mathbb{N}$) shows that
the map%
\begin{align*}
\mathbb{N}^{M}  &  \rightarrow\mathbb{N}^{M-1}\times\mathbb{N},\\
\left(  s_{1},s_{2},\ldots,s_{M}\right)   &  \mapsto\left(  \left(
s_{1},s_{2},\ldots,s_{M-1}\right)  ,s_{M}\right)
\end{align*}
is a bijection. Hence, we can substitute $\left(  \left(  s_{1},s_{2}%
,\ldots,s_{M-1}\right)  ,s_{M}\right)  $ for $\left(  \left(  a_{1}%
,a_{2},\ldots,a_{M-1}\right)  ,q\right)  $ in the sum $\sum_{\substack{\left(
\left(  a_{1},a_{2},\ldots,a_{M-1}\right)  ,q\right)  \in\mathbb{N}%
^{M-1}\times\mathbb{N};\\a_{1}+a_{2}+\cdots+a_{M-1}=r;\\q=n-r}}1$. Thus, we
obtain%
\begin{equation}
\sum_{\substack{\left(  \left(  a_{1},a_{2},\ldots,a_{M-1}\right)  ,q\right)
\in\mathbb{N}^{M-1}\times\mathbb{N};\\a_{1}+a_{2}+\cdots+a_{M-1}%
=r;\\q=n-r}}1=\sum_{\substack{\left(  s_{1},s_{2},\ldots,s_{M}\right)
\in\mathbb{N}^{M};\\s_{1}+s_{2}+\cdots+s_{M-1}=r;\\s_{M}=n-r}}1=\sum
_{\substack{\left(  a_{1},a_{2},\ldots,a_{M}\right)  \in\mathbb{N}^{M}%
;\\a_{1}+a_{2}+\cdots+a_{M-1}=r;\\a_{M}=n-r}}1\nonumber
\end{equation}
(here, we have renamed the summation index $\left(  s_{1},s_{2},\ldots
,s_{M}\right)  $ as $\left(  a_{1},a_{2},\ldots,a_{M}\right)  $). Hence,%
\begin{align*}
\sum_{\substack{\left(  a_{1},a_{2},\ldots,a_{M}\right)  \in\mathbb{N}%
^{M};\\a_{1}+a_{2}+\cdots+a_{M-1}=r;\\a_{M}=n-r}}1  &  =\underbrace{\sum
_{\substack{\left(  \left(  a_{1},a_{2},\ldots,a_{M-1}\right)  ,q\right)
\in\mathbb{N}^{M-1}\times\mathbb{N};\\a_{1}+a_{2}+\cdots+a_{M-1}=r;\\q=n-r}%
}}_{=\sum_{\substack{\left(  a_{1},a_{2},\ldots,a_{M-1}\right)  \in
\mathbb{N}^{M-1};\\a_{1}+a_{2}+\cdots+a_{M-1}=r}}\sum_{\substack{q\in
\mathbb{N};\\q=n-r}}}1=\sum_{\substack{\left(  a_{1},a_{2},\ldots
,a_{M-1}\right)  \in\mathbb{N}^{M-1};\\a_{1}+a_{2}+\cdots+a_{M-1}%
=r}}\underbrace{\sum_{\substack{q\in\mathbb{N};\\q=n-r}}1}%
_{\substack{=1\\\text{(by (\ref{pf.lem.sol.multichoose.sum.sums4}))}}}\\
&  =\sum_{\substack{\left(  a_{1},a_{2},\ldots,a_{M-1}\right)  \in
\mathbb{N}^{M-1};\\a_{1}+a_{2}+\cdots+a_{M-1}=r}}1=\dbinom{r+\left(
M-1\right)  -1}{r}\\
&  \ \ \ \ \ \ \ \ \ \ \left(  \text{by (\ref{pf.lem.sol.multichoose.sum.IH})
(applied to }r\text{ instead of }n\text{)}\right) \\
&  =\dbinom{r+M-2}{r}%
\end{align*}
(since $\left(  M-1\right)  -1=M-2$). Now,%
\[
\underbrace{\sum_{\substack{\left(  a_{1},a_{2},\ldots,a_{M}\right)
\in\mathbb{N}^{M};\\a_{1}+a_{2}+\cdots+a_{M}=n;\\n-a_{M}=r}}}_{\substack{=\sum
_{\substack{\left(  a_{1},a_{2},\ldots,a_{M}\right)  \in\mathbb{N}^{M}%
;\\a_{1}+a_{2}+\cdots+a_{M-1}=r;\\a_{M}=n-r}}\\\text{(by
(\ref{pf.lem.sol.multichoose.sum.sums3}))}}}1=\sum_{\substack{\left(
a_{1},a_{2},\ldots,a_{M}\right)  \in\mathbb{N}^{M};\\a_{1}+a_{2}%
+\cdots+a_{M-1}=r;\\a_{M}=n-r}}1=\dbinom{r+M-2}{r}.
\]
Thus, (\ref{pf.lem.sol.multichoose.sum.sums2}) is proven.]

Now,
\begin{align*}
\underbrace{\sum_{\substack{\left(  a_{1},a_{2},\ldots,a_{M}\right)
\in\mathbb{N}^{M};\\a_{1}+a_{2}+\cdots+a_{M}=n}}}_{\substack{=\sum_{r=0}%
^{n}\sum_{\substack{\left(  a_{1},a_{2},\ldots,a_{M}\right)  \in\mathbb{N}%
^{M};\\a_{1}+a_{2}+\cdots+a_{M}=n;\\n-a_{M}=r}}\\\text{(by
(\ref{pf.lem.sol.multichoose.sum.sums1}))}}}1  &  =\sum_{r=0}^{n}%
\underbrace{\sum_{\substack{\left(  a_{1},a_{2},\ldots,a_{M}\right)
\in\mathbb{N}^{M};\\a_{1}+a_{2}+\cdots+a_{M}=n;\\n-a_{M}=r}}1}%
_{\substack{=\dbinom{r+M-2}{r}\\\text{(by
(\ref{pf.lem.sol.multichoose.sum.sums2}))}}}=\sum_{r=0}^{n}\dbinom{r+M-2}{r}\\
&  =\dbinom{n+\left(  M-2\right)  +1}{n}\\
&  \ \ \ \ \ \ \ \ \ \ \left(  \text{by Lemma \ref{lem.sol.multichoose.hock}
(applied to }q=M-2\text{)}\right) \\
&  =\dbinom{n+M-1}{n}\ \ \ \ \ \ \ \ \ \ \left(  \text{since }\left(
M-2\right)  +1=M-1\right)  .
\end{align*}
In other words, (\ref{pf.lem.sol.multichoose.sum.IG}) holds.

Now, forget that we fixed $n$. We thus have proven that every $n\in\mathbb{N}$
satisfies (\ref{pf.lem.sol.multichoose.sum.IG}). In other words, Lemma
\ref{lem.sol.multichoose.sum} holds for $m=M$. This completes the induction
step. The induction proof of Lemma \ref{lem.sol.multichoose.sum} is thus complete.
\end{proof}

\begin{proof}
[Solution to Exercise \ref{exe.multichoose}.]Lemma
\ref{lem.sol.multichoose.sum} (applied to $m=k$) yields%
\[
\sum_{\substack{\left(  a_{1},a_{2},\ldots,a_{k}\right)  \in\mathbb{N}%
^{k};\\a_{1}+a_{2}+\cdots+a_{k}=n}}1=\dbinom{n+k-1}{n}.
\]
Comparing this with%
\begin{align*}
&  \sum_{\substack{\left(  a_{1},a_{2},\ldots,a_{k}\right)  \in\mathbb{N}%
^{k};\\a_{1}+a_{2}+\cdots+a_{k}=n}}1\\
&  =\left\vert \left\{  \left(  a_{1},a_{2},\ldots,a_{k}\right)  \in
\mathbb{N}^{k}\ \mid\ a_{1}+a_{2}+\cdots+a_{k}=n\right\}  \right\vert \cdot1\\
&  =\left\vert \left\{  \left(  a_{1},a_{2},\ldots,a_{k}\right)  \in
\mathbb{N}^{k}\ \mid\ a_{1}+a_{2}+\cdots+a_{k}=n\right\}  \right\vert \\
&  =\left(  \text{the number of all }k\text{-tuples }\left(  a_{1}%
,a_{2},\ldots,a_{k}\right)  \in\mathbb{N}^{k}\ \text{satisfying}\ a_{1}%
+a_{2}+\cdots+a_{k}=n\right)  ,
\end{align*}
we obtain%
\begin{align*}
&  \left(  \text{the number of all }k\text{-tuples }\left(  a_{1},a_{2}%
,\ldots,a_{k}\right)  \in\mathbb{N}^{k}\ \text{satisfying}\ a_{1}+a_{2}%
+\cdots+a_{k}=n\right) \\
&  =\dbinom{n+k-1}{n}.
\end{align*}
In other words, the number of all $k$-tuples $\left(  a_{1},a_{2},\ldots
,a_{k}\right)  \in\mathbb{N}^{k}$ satisfying $a_{1}+a_{2}+\cdots+a_{k}=n$
equals $\dbinom{n+k-1}{n}$. This solves Exercise \ref{exe.multichoose}.
\end{proof}

\subsection{\label{sect.sol.schmitt-iha-eq9.2}Solution to Exercise
\ref{exe.schmitt-iha-eq9.2}}

Exercise \ref{exe.schmitt-iha-eq9.2} is an identity that tends to creep up in
various seemingly unrelated situations in mathematics. I have first
encountered it in \cite[proof of Theorem 9.5]{Schmitt04} (where it appears
with an incorrect power of $-1$ on the right hand side). It also has recently
appeared on math.stackexchange (\cite{dilemi17}, with $a$, $n-a$ and $k-a$
renamed as $p$, $n$ and $k$), where it has been proven in three different
ways: once using the beta function, once using residues, and once (by myself
in the comments) using finite differences. Let me here give a different,
elementary proof.

We begin with the following identities:

\begin{proposition}
\label{prop.schmitt-iha-eq9.2.binoms}Let $i\in\mathbb{Z}$, $n\in\mathbb{Z}$
and $j\in\mathbb{N}$. Then:

\textbf{(a)} We have
\[
\sum_{k=0}^{j}\left(  -1\right)  ^{k}\dbinom{n}{j-k}\dbinom{k+i-1}{k}%
=\dbinom{n-i}{j}.
\]

\textbf{(b)} If $i$ is positive, then%
\[
\sum_{k=0}^{j}\dfrac{\left(  -1\right)  ^{k}}{k+i}\dbinom{n}{j-k}\dbinom
{k+i}{i}=\dfrac{1}{i}\dbinom{n-i}{j}.
\]

\end{proposition}

\begin{proof}
[Proof of Proposition \ref{prop.schmitt-iha-eq9.2.binoms}.]\textbf{(a)}
Proposition \ref{prop.vandermonde.consequences} \textbf{(a)} (applied to $-i$,
$n$ and $j$ instead of $x$, $y$ and $n$) yields%
\begin{align*}
\dbinom{\left(  -i\right)  +n}{j}  &  =\sum_{k=0}^{j}\underbrace{\dbinom
{-i}{k}}_{\substack{=\left(  -1\right)  ^{k}\dbinom{k-\left(  -i\right)
-1}{k}\\\text{(by Proposition \ref{prop.binom.upper-neg},}\\\text{applied to
}-i\text{ and }k\text{ instead of }m\text{ and }n\text{)}}}\dbinom{n}{j-k}\\
&  =\sum_{k=0}^{j}\left(  -1\right)  ^{k}\underbrace{\dbinom{k-\left(
-i\right)  -1}{k}}_{\substack{=\dbinom{k+i-1}{k}\\\text{(since }k-\left(
-i\right)  -1=k+i-1\text{)}}}\dbinom{n}{j-k}\\
&  =\sum_{k=0}^{j}\left(  -1\right)  ^{k}\dbinom{k+i-1}{k}\dbinom{n}{j-k}%
=\sum_{k=0}^{j}\left(  -1\right)  ^{k}\dbinom{n}{j-k}\dbinom{k+i-1}{k}.
\end{align*}
Thus,%
\[
\sum_{k=0}^{j}\left(  -1\right)  ^{k}\dbinom{n}{j-k}\dbinom{k+i-1}{k}%
=\dbinom{\left(  -i\right)  +n}{j}=\dbinom{n-i}{j}.
\]
This proves Proposition \ref{prop.schmitt-iha-eq9.2.binoms} \textbf{(a)}.

\textbf{(b)} Assume that $i$ is positive. Let $k\in\mathbb{N}$. Then,
$i-1\in\mathbb{N}$ (since $i$ is a positive integer). Thus, $i-1\geq0$. Also,
$\underbrace{k}_{\geq0}+i-1\geq i-1$. Hence, Proposition \ref{prop.binom.symm}
(applied to $k+i-1$ and $i-1$ instead of $m$ and $n$) yields
\[
\dbinom{k+i-1}{i-1}=\dbinom{k+i-1}{\left(  k+i-1\right)  -\left(  i-1\right)
}=\dbinom{k+i-1}{k}%
\]
(since $\left(  k+i-1\right)  -\left(  i-1\right)  =k$).

Furthermore, $i\in\left\{  1,2,3,\ldots\right\}  $ (since $i$ is a positive
integer). Thus, Proposition \ref{prop.binom.X-1} (applied to $k+i$ and $i$
instead of $m$ and $n$) yields $\dbinom{k+i}{i}=\dfrac{k+i}{i}\dbinom
{k+i-1}{i-1}$. Multiplying both sides of this equality by $i$, we find
\[
i\dbinom{k+i}{i}= i\cdot\dfrac{k+i}{i}\dbinom{k+i-1}{i-1} = \left(
k+i\right)  \underbrace{\dbinom{k+i-1}{i-1}}_{=\dbinom{k+i-1}{k}}=\left(
k+i\right)  \dbinom{k+i-1}{k}.
\]
Hence,
\[
\dbinom{k+i}{i}=\dfrac{1}{i}\left(  k+i\right)  \dbinom{k+i-1}{k}.
\]
Thus,%
\begin{align}
&  \dfrac{\left(  -1\right)  ^{k}}{k+i}\dbinom{n}{j-k}\underbrace{\dbinom
{k+i}{i}}_{=\dfrac{1}{i}\left(  k+i\right)  \dbinom{k+i-1}{k}}\nonumber\\
&  =\dfrac{\left(  -1\right)  ^{k}}{k+i}\dbinom{n}{j-k}\cdot\dfrac{1}%
{i}\left(  k+i\right)  \dbinom{k+i-1}{k}\nonumber\\
&  =\dfrac{1}{i}\left(  -1\right)  ^{k}\dbinom{n}{j-k}\dbinom{k+i-1}{k}.
\label{pf.prop.schmitt-iha-eq9.2.binoms.b.3}%
\end{align}

Now, forget that we fixed $k$. We thus have proven
(\ref{pf.prop.schmitt-iha-eq9.2.binoms.b.3}) for each $k\in\mathbb{N}$. Now,%
\begin{align*}
&  \sum_{k=0}^{j}\underbrace{\dfrac{\left(  -1\right)  ^{k}}{k+i}\dbinom
{n}{j-k}\dbinom{k+i}{i}}_{\substack{=\dfrac{1}{i}\left(  -1\right)
^{k}\dbinom{n}{j-k}\dbinom{k+i-1}{k}\\\text{(by
(\ref{pf.prop.schmitt-iha-eq9.2.binoms.b.3}))}}}\\
&  =\sum_{k=0}^{j}\dfrac{1}{i}\left(  -1\right)  ^{k}\dbinom{n}{j-k}%
\dbinom{k+i-1}{k}\\
&  =\dfrac{1}{i}\underbrace{\sum_{k=0}^{j}\left(  -1\right)  ^{k}\dbinom
{n}{j-k}\dbinom{k+i-1}{k}}_{\substack{=\dbinom{n-i}{j}\\\text{(by Proposition
\ref{prop.schmitt-iha-eq9.2.binoms} \textbf{(a)})}}}=\dfrac{1}{i}\dbinom
{n-i}{j}.
\end{align*}
This proves Proposition \ref{prop.schmitt-iha-eq9.2.binoms} \textbf{(b)}.
\end{proof}

Let us now solve the actual exercise:

\begin{proof}
[Solution to Exercise \ref{exe.schmitt-iha-eq9.2}.]From $n\geq a$, we obtain
$n-a\in\mathbb{N}$. Also, $n\geq a\geq1\geq0$, and therefore Proposition
\ref{prop.binom.formula} (applied to $n$ and $a$ instead of $m$ and $n$)
yields $\dbinom{n}{a}=\dfrac{n!}{a!\left(  n-a\right)  !}\neq0$ (since
$n!\neq0$).

Any $k\in\left\{  a,a+1,\ldots,n\right\}  $ satisfies%
\begin{equation}
\dbinom{n}{a}\dbinom{n-a}{k-a}=\dbinom{n}{n-k}\dbinom{k}{a}
\label{sol.schmitt-iha-eq9.2.1}%
\end{equation}
\footnote{\textit{Proof of (\ref{sol.schmitt-iha-eq9.2.1}):} Let $k\in\left\{
a,a+1,\ldots,n\right\}  $. Thus, $a\leq k\leq n$, so that $k\geq a\geq1\geq0$,
so that $k\in\mathbb{N}$. Hence, Proposition \ref{prop.binom.symm} (applied to
$n$ and $k$ instead of $m$ and $n$) yields $\dbinom{n}{k}=\dbinom{n}{n-k}$
(since $n\geq k$).
\par
But Proposition \ref{prop.binom.trinom-rev} (applied to $n$ and $k$ instead of
$m$ and $i$) shows that%
\[
\dbinom{n}{k}\dbinom{k}{a}=\dbinom{n}{a}\dbinom{n-a}{k-a}.
\]
Hence,%
\[
\dbinom{n}{a}\dbinom{n-a}{k-a}=\underbrace{\dbinom{n}{k}}_{=\dbinom{n}{n-k}%
}\dbinom{k}{a}=\dbinom{n}{n-k}\dbinom{k}{a}.
\]
This proves (\ref{sol.schmitt-iha-eq9.2.1}).}.

We have%
\begin{align*}
\dbinom{n}{a}\sum_{k=a}^{n}\dfrac{\left(  -1\right)  ^{k}}{k}\dbinom
{n-a}{k-a}  &  =\sum_{k=a}^{n}\dfrac{\left(  -1\right)  ^{k}}{k}%
\underbrace{\dbinom{n}{a}\dbinom{n-a}{k-a}}_{\substack{=\dbinom{n}{n-k}%
\dbinom{k}{a}\\\text{(by (\ref{sol.schmitt-iha-eq9.2.1}))}}}=\sum_{k=a}%
^{n}\dfrac{\left(  -1\right)  ^{k}}{k}\dbinom{n}{n-k}\dbinom{k}{a}\\
&  =\sum_{k=0}^{n-a}\underbrace{\dfrac{\left(  -1\right)  ^{k+a}}{k+a}%
}_{\substack{=\dfrac{\left(  -1\right)  ^{k}\left(  -1\right)  ^{a}}%
{k+a}\\\text{(since }\left(  -1\right)  ^{k+a}=\left(  -1\right)  ^{k}\left(
-1\right)  ^{a}\text{)}}}\underbrace{\dbinom{n}{n-\left(  k+a\right)  }%
}_{\substack{=\dbinom{n}{\left(  n-a\right)  -k}\\\text{(since }n-\left(
k+a\right)  =\left(  n-a\right)  -k\text{)}}}\dbinom{k+a}{a}\\
&  \ \ \ \ \ \ \ \ \ \ \left(  \text{here, we have substituted }k+a\text{ for
}k\text{ in the sum}\right) \\
&  =\sum_{k=0}^{n-a}\dfrac{\left(  -1\right)  ^{k}\left(  -1\right)  ^{a}%
}{k+a}\dbinom{n}{\left(  n-a\right)  -k}\dbinom{k+a}{a}\\
&  =\left(  -1\right)  ^{a}\underbrace{\sum_{k=0}^{n-a}\dfrac{\left(
-1\right)  ^{k}}{k+a}\dbinom{n}{\left(  n-a\right)  -k}\dbinom{k+a}{a}%
}_{\substack{=\dfrac{1}{a}\dbinom{n-a}{n-a}\\\text{(by Proposition
\ref{prop.schmitt-iha-eq9.2.binoms} \textbf{(b)},}\\\text{applied to
}j=n-a\text{ and }i=a\text{)}}}\\
&  =\left(  -1\right)  ^{a}\dfrac{1}{a}\underbrace{\dbinom{n-a}{n-a}%
}_{\substack{=1\\\text{(by Proposition \ref{prop.binom.mm}}\\\text{(applied to
}m=n-a\text{))}}}=\left(  -1\right)  ^{a}\dfrac{1}{a}.
\end{align*}
We can divide both sides of this equality by $\dbinom{n}{a}$ (since
$\dbinom{n}{a}\neq0$). Thus, we find%
\[
\sum_{k=a}^{n}\dfrac{\left(  -1\right)  ^{k}}{k}\dbinom{n-a}{k-a}=\left(
-1\right)  ^{a}\dfrac{1}{a}/\dbinom{n}{a}=\dfrac{\left(  -1\right)  ^{a}%
}{a\dbinom{n}{a}}.
\]
This solves Exercise \ref{exe.schmitt-iha-eq9.2}.
\end{proof}

\subsection{\label{sect.sol.bininv}Solution to Exercise \ref{exe.bininv}}

We shall now prepare for the solution of Exercise \ref{exe.bininv}.

Let us first fix some notations.

\begin{definition}
Let $N\in\mathbb{N}$. We shall consider $N$ to be fixed for the whole Section
\ref{sect.sol.bininv}.

Throughout Section \ref{sect.sol.bininv}, we shall use the word
\textquotedblleft list\textquotedblright\ for an $\left(  N+1\right)  $-tuple
of rational numbers. In other words, a \textquotedblleft
list\textquotedblright\ will mean an element of $\mathbb{Q}^{N+1}$. Thus, for
example, when we say \textquotedblleft the list $\left(  1,1,\ldots,1\right)
$\textquotedblright, we mean the list $\left(  \underbrace{1,1,\ldots
,1}_{N+1\text{ entries}}\right)  $ (because any list has to be an $\left(
N+1\right)  $-tuple).
\end{definition}

\begin{definition}
\label{def.sol.bininv.bintr}The \textit{binomial transform} of a list $\left(
f_{0},f_{1},\ldots,f_{N}\right)  \in\mathbb{Q}^{N+1}$ is defined to be the
list $\left(  g_{0},g_{1},\ldots,g_{N}\right)  $ defined by%
\[
\left(  g_{n}=\sum_{i=0}^{n}\left(  -1\right)  ^{i}\dbinom{n}{i}%
f_{i}\ \ \ \ \ \ \ \ \ \ \text{for every }n\in\left\{  0,1,\ldots,N\right\}
\right)  .
\]

\end{definition}

Clearly, Definition \ref{def.sol.bininv.bintr} generalizes the definition of
the binomial transform we gave in Exercise \ref{exe.bininv} (because any
finite sequence $\left(  f_{0},f_{1},\ldots,f_{N}\right)  \in\mathbb{Z}^{N+1}$
of integers is clearly a list in $\mathbb{Q}^{N+1}$).

We shall use the Iverson bracket notation introduced in Definition
\ref{def.iverson}.

The following fact is easy:

\begin{proposition}
\label{prop.binom.1-1}Let $m\in\mathbb{N}$. Then,%
\[
\sum_{k=0}^{m}\left(  -1\right)  ^{k}\dbinom{m}{k}=\left[  m=0\right]  .
\]

\end{proposition}

\begin{vershort}
\begin{proof}
[Proof of Proposition \ref{prop.binom.1-1}.]Proposition
\ref{prop.binom.bin-id} \textbf{(c)} (applied to $n=m$) yields%
\[
\sum_{k=0}^{m}\left(  -1\right)  ^{k}\dbinom{m}{k}=%
\begin{cases}
1, & \text{if }m=0;\\
0, & \text{if }m\neq0
\end{cases}
=%
\begin{cases}
1, & \text{if }m=0\text{ is true};\\
0, & \text{if }m=0\text{ is false}%
\end{cases}
.
\]
Comparing this with%
\[
\left[  m=0\right]  =%
\begin{cases}
1, & \text{if }m=0\text{ is true};\\
0, & \text{if }m=0\text{ is false}%
\end{cases}
\ \ \ \ \ \ \ \ \ \ \left(  \text{by the definition of }\left[  m=0\right]
\right)  ,
\]
we obtain $\sum_{k=0}^{m}\left(  -1\right)  ^{k}\dbinom{m}{k}=\left[
m=0\right]  $. This proves Proposition \ref{prop.binom.1-1}.
\end{proof}
\end{vershort}

\begin{verlong}
\begin{proof}
[Proof of Proposition \ref{prop.binom.1-1}.]Proposition
\ref{prop.binom.bin-id} \textbf{(c)} (applied to $n=m$) yields%
\[
\sum_{k=0}^{m}\left(  -1\right)  ^{k}\dbinom{m}{k}=%
\begin{cases}
1, & \text{if }m=0;\\
0, & \text{if }m\neq0
\end{cases}
.
\]
Comparing this with%
\begin{align*}
\left[  m=0\right]   &  =%
\begin{cases}
1, & \text{if }m=0\text{ is true};\\
0, & \text{if }m=0\text{ is false}%
\end{cases}
\ \ \ \ \ \ \ \ \ \ \left(  \text{by the definition of }\left[  m=0\right]
\right) \\
&  =%
\begin{cases}
1, & \text{if }m=0;\\
0, & \text{if }m\neq0
\end{cases}
\\
&  \ \ \ \ \ \ \ \ \ \ \left(
\begin{array}
[c]{c}%
\text{since \textquotedblleft}m=0\text{ is true\textquotedblright\ means the
same as \textquotedblleft}m=0\text{\textquotedblright,}\\
\text{whereas \textquotedblleft}m=0\text{ is false\textquotedblright\ means
the same as \textquotedblleft}m\neq0\text{\textquotedblright}%
\end{array}
\right)  ,
\end{align*}
we obtain $\sum_{k=0}^{m}\left(  -1\right)  ^{k}\dbinom{m}{k}=\left[
m=0\right]  $. This proves Proposition \ref{prop.binom.1-1}.
\end{proof}
\end{verlong}

\begin{corollary}
\label{cor.binom.1-1cor}Let $n\in\mathbb{N}$. Let $i\in\left\{  0,1,\ldots
,n\right\}  $. Then,%
\[
\sum_{j=i}^{n}\left(  -1\right)  ^{j+i}\dbinom{n}{j}\dbinom{j}{i}=\left[
i=n\right]  .
\]

\end{corollary}

\begin{proof}
[Proof of Corollary \ref{cor.binom.1-1cor}.]We have $i\in\left\{
0,1,\ldots,n\right\}  $. Hence, $i\leq n$. Also, $i\in\left\{  0,1,\ldots
,n\right\}  \subseteq\mathbb{N}$. From $i\leq n$, we obtain $n-i\geq0$. Thus,
$n-i\in\mathbb{N}$. Therefore, Proposition \ref{prop.binom.1-1} (applied to
$m=n-i$) yields%
\begin{equation}
\sum_{k=0}^{n-i}\left(  -1\right)  ^{k}\dbinom{n-i}{k}=\left[  n-i=0\right]  .
\label{pf.cor.binom.1-1cor.0}%
\end{equation}

Let $j\in\left\{  i,i+1,\ldots,n\right\}  $. Then, $j\geq i$, so that $j\geq
i\geq0$ and thus $j\in\mathbb{N}$. Hence, Proposition
\ref{prop.binom.trinom-rev} (applied to $n$, $j$ and $i$ instead of $m$, $i$
and $a$) yields%
\begin{equation}
\dbinom{n}{j}\dbinom{j}{i}=\dbinom{n}{i}\dbinom{n-i}{j-i}.
\label{pf.cor.binom.1-1cor.1}%
\end{equation}
Also, $j+i\equiv j-i\operatorname{mod}2$ (since $\left(  j+i\right)  -\left(
j-i\right)  =2i$ is even). Thus, $\left(  -1\right)  ^{j+i}=\left(  -1\right)
^{j-i}$. Multiplying this equality with the equality
(\ref{pf.cor.binom.1-1cor.1}), we obtain%
\begin{equation}
\left(  -1\right)  ^{j+i}\dbinom{n}{j}\dbinom{j}{i}=\left(  -1\right)
^{j-i}\dbinom{n}{i}\dbinom{n-i}{j-i}. \label{pf.cor.binom.1-1cor.2}%
\end{equation}

Now, forget that we fixed $j$. We thus have proven
(\ref{pf.cor.binom.1-1cor.2}) for each $j\in\left\{  i,i+1,\ldots,n\right\}
$. Hence,%
\begin{align}
&  \sum_{j=i}^{n}\underbrace{\left(  -1\right)  ^{j+i}\dbinom{n}{j}\dbinom
{j}{i}}_{\substack{=\left(  -1\right)  ^{j-i}\dbinom{n}{i}\dbinom{n-i}%
{j-i}\\\text{(by (\ref{pf.cor.binom.1-1cor.2}))}}}\nonumber\\
&  =\sum_{j=i}^{n}\left(  -1\right)  ^{j-i}\dbinom{n}{i}\dbinom{n-i}{j-i}%
=\sum_{k=0}^{n-i}\left(  -1\right)  ^{k}\dbinom{n}{i}\dbinom{n-i}%
{k}\nonumber\\
&  \ \ \ \ \ \ \ \ \ \ \left(  \text{here, we have substituted }k\text{ for
}j-i\text{ in the sum}\right) \nonumber\\
&  =\dbinom{n}{i}\underbrace{\sum_{k=0}^{n-i}\left(  -1\right)  ^{k}%
\dbinom{n-i}{k}}_{\substack{=\left[  n-i=0\right]  \\\text{(by
(\ref{pf.cor.binom.1-1cor.0}))}}}=\dbinom{n}{i}\left[  n-i=0\right]  .
\label{pf.cor.binom.1-1cor.4}%
\end{align}
But it is easy to see that%
\begin{equation}
\dbinom{n}{i}\left[  n-i=0\right]  =\left[  i=n\right]
\label{pf.cor.binom.1-1cor.5}%
\end{equation}
\footnote{\textit{Proof of (\ref{pf.cor.binom.1-1cor.5}):} We are in one of
the following two cases:
\par
\textit{Case 1:} We have $i\neq n$.
\par
\textit{Case 2:} We have $i=n$.
\par
Let us consider Case 1 first. In this case, we have $i\neq n$. In other words,
$n\neq i$. Hence, $n-i\neq0$. Thus, the statement $n-i=0$ is false. Hence,
$\left[  n-i=0\right]  =0$, so that $\dbinom{n}{i}\underbrace{\left[
n-i=0\right]  }_{=0}=0$. Comparing this with $\left[  i=n\right]  =0$ (since
we don't have $i=n$ (since $i\neq n$)), we obtain $\dbinom{n}{i}\left[
n-i=0\right]  =\left[  i=n\right]  $. Hence, (\ref{pf.cor.binom.1-1cor.5}) is
proven in Case 1.
\par
Now, let us consider Case 2. In this case, we have $i=n$. In other words,
$n=i$. Thus, $n-i=0$. Thus, $\left[  n-i=0\right]  =1$. Also, from $n=i$, we
obtain $\dbinom{n}{i}=\dbinom{i}{i}=1$ (by Proposition \ref{prop.binom.mm}
(applied to $m=i$)). Hence, $\underbrace{\dbinom{n}{i}}_{=1}%
\underbrace{\left[  n-i=0\right]  }_{=1}=1$. Comparing this with $\left[
i=n\right]  =1$ (since $i=n$), we obtain $\dbinom{n}{i}\left[  n-i=0\right]
=\left[  i=n\right]  $. Hence, (\ref{pf.cor.binom.1-1cor.5}) is proven in Case
2.
\par
We thus have proven (\ref{pf.cor.binom.1-1cor.5}) in each of the two Cases 1
and 2. Thus, (\ref{pf.cor.binom.1-1cor.5}) always holds.}. Thus,
(\ref{pf.cor.binom.1-1cor.4}) becomes%
\[
\sum_{j=i}^{n}\left(  -1\right)  ^{j+i}\dbinom{n}{j}\dbinom{j}{i}=\dbinom
{n}{i}\left[  n-i=0\right]  =\left[  i=n\right]  .
\]
This proves Corollary \ref{cor.binom.1-1cor}.
\end{proof}

\begin{proposition}
\label{prop.sol.bintra.a0}Let $\left(  a_{0},a_{1},\ldots,a_{N}\right)  $ and
$\left(  b_{0},b_{1},\ldots,b_{N}\right)  $ be two $\left(  N+1\right)
$-tuples of rational numbers. Assume that
\[
b_{n}=\sum_{i=0}^{n}\left(  -1\right)  ^{i}\dbinom{n}{i}a_{i}%
\ \ \ \ \ \ \ \ \ \ \text{for each }n\in\left\{  0,1,\ldots,N\right\}  .
\]
Then,
\[
a_{n}=\sum_{i=0}^{n}\left(  -1\right)  ^{i}\dbinom{n}{i}b_{i}%
\ \ \ \ \ \ \ \ \ \ \text{for each }n\in\left\{  0,1,\ldots,N\right\}  .
\]

\end{proposition}

\begin{vershort}
\begin{proof}
[Proof of Proposition \ref{prop.sol.bintra.a0}.]We have assumed that%
\begin{equation}
b_{n}=\sum_{i=0}^{n}\left(  -1\right)  ^{i}\dbinom{n}{i}a_{i}
\label{pf.prop.sol.bintra.a0.short.1}%
\end{equation}
for each $n\in\left\{  0,1,\ldots,N\right\}  $.

Now, let $n\in\left\{  0,1,\ldots,N\right\}  $. Then,
\begin{align*}
&  \sum_{i=0}^{n}\left(  -1\right)  ^{i}\dbinom{n}{i}b_{i}\\
&  =\sum_{j=0}^{n}\left(  -1\right)  ^{j}\dbinom{n}{j}\underbrace{b_{j}%
}_{\substack{=\sum_{i=0}^{j}\left(  -1\right)  ^{i}\dbinom{j}{i}%
a_{i}\\\text{(by (\ref{pf.prop.sol.bintra.a0.short.1}) (applied to
}j\\\text{instead of }n\text{))}}}\ \ \ \ \ \ \ \ \ \ \left(
\begin{array}
[c]{c}%
\text{here, we have renamed the}\\
\text{summation index }i\text{ as }j
\end{array}
\right) \\
&  =\sum_{j=0}^{n}\left(  -1\right)  ^{j}\dbinom{n}{j}\sum_{i=0}^{j}\left(
-1\right)  ^{i}\dbinom{j}{i}a_{i}\\
&  =\underbrace{\sum_{j=0}^{n}}_{=\sum_{j\in\left\{  0,1,\ldots,n\right\}  }%
}\underbrace{\sum_{i=0}^{j}}_{\substack{=\sum_{i\in\left\{  0,1,\ldots
,j\right\}  }=\sum_{\substack{i\in\left\{  0,1,\ldots,n\right\}  ;\\i\leq
j}}\\\text{(since the elements of }\left\{  0,1,\ldots,j\right\}  \\\text{are
precisely the elements }i\in\left\{  0,1,\ldots,n\right\}  \\\text{satisfying
}i\leq j\text{ (because }j\leq n\text{))}}}\left(  -1\right)  ^{j}%
\underbrace{\dbinom{n}{j}\left(  -1\right)  ^{i}}_{=\left(  -1\right)
^{i}\dbinom{n}{j}}\dbinom{j}{i}a_{i}\\
&  =\underbrace{\sum_{j\in\left\{  0,1,\ldots,n\right\}  }\sum_{\substack{i\in
\left\{  0,1,\ldots,n\right\}  ;\\i\leq j}}}_{=\sum_{i\in\left\{
0,1,\ldots,n\right\}  }\sum_{\substack{j\in\left\{  0,1,\ldots,n\right\}
;\\i\leq j}}}\underbrace{\left(  -1\right)  ^{j}\left(  -1\right)  ^{i}%
}_{=\left(  -1\right)  ^{j+i}}\dbinom{n}{j}\dbinom{j}{i}a_{i}\\
&  =\sum_{i\in\left\{  0,1,\ldots,n\right\}  }\underbrace{\sum_{\substack{j\in
\left\{  0,1,\ldots,n\right\}  ;\\i\leq j}}}_{=\sum_{j=i}^{n}}\left(
-1\right)  ^{j+i}\dbinom{n}{j}\dbinom{j}{i}a_{i}\\
&  =\sum_{i\in\left\{  0,1,\ldots,n\right\}  }\sum_{j=i}^{n}\left(  -1\right)
^{j+i}\dbinom{n}{j}\dbinom{j}{i}a_{i}=\sum_{i\in\left\{  0,1,\ldots,n\right\}
}\underbrace{\left(  \sum_{j=i}^{n}\left(  -1\right)  ^{j+i}\dbinom{n}%
{j}\dbinom{j}{i}\right)  }_{\substack{=\left[  i=n\right]  \\\text{(by
Corollary \ref{cor.binom.1-1cor})}}}a_{i}\\
&  =\sum_{i\in\left\{  0,1,\ldots,n\right\}  }\left[  i=n\right]
a_{i}=\underbrace{\left[  n=n\right]  }_{\substack{=1\\\text{(since
}n=n\text{)}}}a_{n}+\sum_{\substack{i\in\left\{  0,1,\ldots,n\right\}
;\\i\neq n}}\underbrace{\left[  i=n\right]  }_{\substack{=0\\\text{(since
}i\neq n\text{)}}}a_{i}\\
&  \ \ \ \ \ \ \ \ \ \ \left(  \text{here, we have split off the addend for
}i=n\text{ from the sum}\right) \\
&  =a_{n}+\underbrace{\sum_{\substack{i\in\left\{  0,1,\ldots,n\right\}
;\\i\neq n}}0a_{i}}_{=0}=a_{n}.
\end{align*}
In other words, $a_{n}=\sum_{i=0}^{n}\left(  -1\right)  ^{i}\dbinom{n}{i}%
b_{i}$. This proves Proposition \ref{prop.sol.bintra.a0}.
\end{proof}
\end{vershort}

\begin{verlong}
\begin{proof}
[Proof of Proposition \ref{prop.sol.bintra.a0}.]We have assumed that%
\begin{equation}
b_{n}=\sum_{i=0}^{n}\left(  -1\right)  ^{i}\dbinom{n}{i}a_{i}
\label{pf.prop.sol.bintra.a0.1}%
\end{equation}
for each $n\in\left\{  0,1,\ldots,N\right\}  $.

Now, let $n\in\left\{  0,1,\ldots,N\right\}  $. Thus, $n\geq0$ and $n\leq N$.

For each $j\in\left\{  0,1,\ldots,n\right\}  $, we have the following equality
of summation signs:%
\begin{equation}
\sum_{i=0}^{j}=\sum_{\substack{i\in\left\{  0,1,\ldots,n\right\}  ;\\i\leq j}}
\label{pf.prop.sol.bintra.a0.2}%
\end{equation}
\footnote{\textit{Proof of (\ref{pf.prop.sol.bintra.a0.2}):} Let $j\in\left\{
0,1,\ldots,n\right\}  $. Thus, $0\leq j\leq n$. From $j\leq n$, we obtain
$\left\{  0,1,\ldots,j\right\}  \subseteq\left\{  0,1,\ldots,n\right\}  $.
Thus, $\left\{  0,1,\ldots,n\right\}  \cap\left\{  0,1,\ldots,j\right\}
=\left\{  0,1,\ldots,j\right\}  $. Now, we have the following equality of
summation signs:%
\begin{align*}
\sum_{i=0}^{j}  &  =\sum_{i\in\left\{  0,1,\ldots,j\right\}  }=\sum
_{i\in\left\{  0,1,\ldots,n\right\}  \cap\left\{  0,1,\ldots,j\right\}
}\ \ \ \ \ \ \ \ \ \ \left(  \text{since }\left\{  0,1,\ldots,j\right\}
=\left\{  0,1,\ldots,n\right\}  \cap\left\{  0,1,\ldots,j\right\}  \right) \\
&  =\sum_{\substack{i\in\left\{  0,1,\ldots,n\right\}  ;\\i\in\left\{
0,1,\ldots,j\right\}  }}=\sum_{\substack{i\in\left\{  0,1,\ldots,n\right\}
;\\i\leq j}}
\end{align*}
(because for an $i\in\left\{  0,1,\ldots,n\right\}  $, the condition $\left(
i\in\left\{  0,1,\ldots,j\right\}  \right)  $ is equivalent to the condition
$\left(  i\leq j\right)  $ (since $i\in\left\{  0,1,\ldots,n\right\}
\subseteq\mathbb{N}$)). This proves (\ref{pf.prop.sol.bintra.a0.2}).}. Also,
for each $j\in\left\{  0,1,\ldots,n\right\}  $, we have%
\begin{equation}
b_{j}=\sum_{\substack{i\in\left\{  0,1,\ldots,n\right\}  ;\\i\leq j}}\left(
-1\right)  ^{i}\dbinom{j}{i}a_{i} \label{pf.prop.sol.bintra.a0.3}%
\end{equation}
\footnote{\textit{Proof of (\ref{pf.prop.sol.bintra.a0.3}):} Let $j\in\left\{
0,1,\ldots,n\right\}  $. Then, $j\in\left\{  0,1,\ldots,n\right\}
\subseteq\left\{  0,1,\ldots,N\right\}  $ (since $n\leq N$). Hence,
(\ref{pf.prop.sol.bintra.a0.1}) (applied to $j$ instead of $n$) yields%
\[
b_{j}=\underbrace{\sum_{i=0}^{j}}_{\substack{=\sum_{\substack{i\in\left\{
0,1,\ldots,n\right\}  ;\\i\leq j}}\\\text{(by (\ref{pf.prop.sol.bintra.a0.2}%
))}}}\left(  -1\right)  ^{i}\dbinom{j}{i}a_{i}=\sum_{\substack{i\in\left\{
0,1,\ldots,n\right\}  ;\\i\leq j}}\left(  -1\right)  ^{i}\dbinom{j}{i}a_{i}.
\]
This proves (\ref{pf.prop.sol.bintra.a0.3}).}.

Now,%
\begin{align*}
&  \sum_{i=0}^{n}\left(  -1\right)  ^{i}\dbinom{n}{i}b_{i}\\
&  =\underbrace{\sum_{j=0}^{n}}_{=\sum_{j\in\left\{  0,1,\ldots,n\right\}  }%
}\left(  -1\right)  ^{j}\dbinom{n}{j}\underbrace{b_{j}}_{\substack{=\sum
_{\substack{i\in\left\{  0,1,\ldots,n\right\}  ;\\i\leq j}}\left(  -1\right)
^{i}\dbinom{j}{i}a_{i}\\\text{(by (\ref{pf.prop.sol.bintra.a0.3}))}}}\\
&  \ \ \ \ \ \ \ \ \ \ \left(  \text{here, we have renamed the summation index
}i\text{ as }j\right) \\
&  =\sum_{j\in\left\{  0,1,\ldots,n\right\}  }\left(  -1\right)  ^{j}%
\dbinom{n}{j}\sum_{\substack{i\in\left\{  0,1,\ldots,n\right\}  ;\\i\leq
j}}\left(  -1\right)  ^{i}\dbinom{j}{i}a_{i}\\
&  =\underbrace{\sum_{j\in\left\{  0,1,\ldots,n\right\}  }\sum_{\substack{i\in
\left\{  0,1,\ldots,n\right\}  ;\\i\leq j}}}_{=\sum_{i\in\left\{
0,1,\ldots,n\right\}  }\sum_{\substack{j\in\left\{  0,1,\ldots,n\right\}
;\\i\leq j}}}\left(  -1\right)  ^{j}\underbrace{\dbinom{n}{j}\left(
-1\right)  ^{i}}_{=\left(  -1\right)  ^{i}\dbinom{n}{j}}\dbinom{j}{i}a_{i}\\
&  =\sum_{i\in\left\{  0,1,\ldots,n\right\}  }\underbrace{\sum_{\substack{j\in
\left\{  0,1,\ldots,n\right\}  ;\\i\leq j}}}_{=\sum_{j=i}^{n}}%
\underbrace{\left(  -1\right)  ^{j}\left(  -1\right)  ^{i}}_{=\left(
-1\right)  ^{j+i}}\dbinom{n}{j}\dbinom{j}{i}a_{i}\\
&  =\sum_{i\in\left\{  0,1,\ldots,n\right\}  }\sum_{j=i}^{n}\left(  -1\right)
^{j+i}\dbinom{n}{j}\dbinom{j}{i}a_{i}=\sum_{i\in\left\{  0,1,\ldots,n\right\}
}\underbrace{\left(  \sum_{j=i}^{n}\left(  -1\right)  ^{j+i}\dbinom{n}%
{j}\dbinom{j}{i}\right)  }_{\substack{=\left[  i=n\right]  \\\text{(by
Corollary \ref{cor.binom.1-1cor})}}}a_{i}\\
&  =\sum_{i\in\left\{  0,1,\ldots,n\right\}  }\left[  i=n\right]
a_{i}=\underbrace{\left[  n=n\right]  }_{\substack{=1\\\text{(since
}n=n\text{)}}}a_{n}+\sum_{\substack{i\in\left\{  0,1,\ldots,n\right\}
;\\i\neq n}}\underbrace{\left[  i=n\right]  }_{\substack{=0\\\text{(since
}i\neq n\text{)}}}a_{i}\\
&  \ \ \ \ \ \ \ \ \ \ \left(
\begin{array}
[c]{c}%
\text{here, we have split off the addend for }i=n\text{ from the sum}\\
\text{(since }n\in\left\{  0,1,\ldots,n\right\}  \text{ (because }%
n\geq0\text{))}%
\end{array}
\right) \\
&  =a_{n}+\underbrace{\sum_{\substack{i\in\left\{  0,1,\ldots,n\right\}
;\\i\neq n}}0a_{i}}_{=0}=a_{n}.
\end{align*}
In other words, $a_{n}=\sum_{i=0}^{n}\left(  -1\right)  ^{i}\dbinom{n}{i}%
b_{i}$. This proves Proposition \ref{prop.sol.bintra.a0}.
\end{proof}
\end{verlong}

From this, it is easy to derive the following statement, which generalizes
Exercise \ref{exe.bininv} \textbf{(a)}:

\begin{corollary}
\label{cor.sol.bininv.a}Let $\left(  f_{0},f_{1},\ldots,f_{N}\right)
\in\mathbb{Q}^{N+1}$ be a list. Let $\left(  g_{0},g_{1},\ldots,g_{N}\right)
$ be the binomial transform of $\left(  f_{0},f_{1},\ldots,f_{N}\right)  $.
Then, $\left(  f_{0},f_{1},\ldots,f_{N}\right)  $ is the binomial transform of
$\left(  g_{0},g_{1},\ldots,g_{N}\right)  $.
\end{corollary}

\begin{proof}
[Proof of Corollary \ref{cor.sol.bininv.a}.]We know that $\left(  g_{0}%
,g_{1},\ldots,g_{N}\right)  $ is the binomial transform of $\left(
f_{0},f_{1},\ldots,f_{N}\right)  $. In other words, we have%
\[
g_{n}=\sum_{i=0}^{n}\left(  -1\right)  ^{i}\dbinom{n}{i}f_{i}%
\ \ \ \ \ \ \ \ \ \ \text{for every }n\in\left\{  0,1,\ldots,N\right\}
\]
(by the definition of the binomial transform of a list). Hence, Proposition
\ref{prop.sol.bintra.a0} (applied to $a_{i}=f_{i}$ and $b_{i}=g_{i}$) yields
that%
\begin{equation}
f_{n}=\sum_{i=0}^{n}\left(  -1\right)  ^{i}\dbinom{n}{i}g_{i}%
\ \ \ \ \ \ \ \ \ \ \text{for each }n\in\left\{  0,1,\ldots,N\right\}  .
\label{pf.cor.sol.bininv.a.1}%
\end{equation}

Let $\left(  h_{0},h_{1},\ldots,h_{N}\right)  $ be the binomial transform of
$\left(  g_{0},g_{1},\ldots,g_{N}\right)  $. Thus,%
\begin{equation}
h_{n}=\sum_{i=0}^{n}\left(  -1\right)  ^{i}\dbinom{n}{i}g_{i}%
\ \ \ \ \ \ \ \ \ \ \text{for every }n\in\left\{  0,1,\ldots,N\right\}
\label{pf.cor.sol.bininv.a.2}%
\end{equation}
(by the definition of the binomial transform of a list). Hence, each
$n\in\left\{  0,1,\ldots,N\right\}  $ satisfies%
\[
h_{n}=\sum_{i=0}^{n}\left(  -1\right)  ^{i}\dbinom{n}{i}g_{i}=f_{n}%
\ \ \ \ \ \ \ \ \ \ \left(  \text{by (\ref{pf.cor.sol.bininv.a.1})}\right)  .
\]
In other words, $\left(  h_{0},h_{1},\ldots,h_{N}\right)  =\left(  f_{0}%
,f_{1},\ldots,f_{N}\right)  $.

Recall that $\left(  h_{0},h_{1},\ldots,h_{N}\right)  $ is the binomial
transform of $\left(  g_{0},g_{1},\ldots,g_{N}\right)  $. In other words,
$\left(  f_{0},f_{1},\ldots,f_{N}\right)  $ is the binomial transform of
$\left(  g_{0},g_{1},\ldots,g_{N}\right)  $ (since \newline$\left(
h_{0},h_{1},\ldots,h_{N}\right)  =\left(  f_{0},f_{1},\ldots,f_{N}\right)  $).
This proves Corollary \ref{cor.sol.bininv.a}.
\end{proof}

Next, let us state a simple consequence of the binomial formula:

\begin{corollary}
\label{cor.binom.(1-q)n}Let $n\in\mathbb{N}$ and $q\in\mathbb{Q}$. Then,%
\[
\sum_{i=0}^{n}\left(  -1\right)  ^{i}\dbinom{n}{i}q^{i}=\left(  1-q\right)
^{n}.
\]

\end{corollary}

\begin{proof}
[Proof of Corollary \ref{cor.binom.(1-q)n}.]Proposition
\ref{prop.binom.binomial} (applied to $x=-q$ and $y=1$) yields
\begin{align*}
\left(  \left(  -q\right)  +1\right)  ^{n}  &  =\sum_{k=0}^{n}\dbinom{n}%
{k}\left(  -q\right)  ^{k}\underbrace{1^{n-k}}_{=1}=\sum_{k=0}^{n}\dbinom
{n}{k}\underbrace{\left(  -q\right)  ^{k}}_{=\left(  -1\right)  ^{k}q^{k}}\\
&  =\sum_{k=0}^{n}\underbrace{\dbinom{n}{k}\left(  -1\right)  ^{k}}_{=\left(
-1\right)  ^{k}\dbinom{n}{k}}q^{k}=\sum_{k=0}^{n}\left(  -1\right)
^{k}\dbinom{n}{k}q^{k}=\sum_{i=0}^{n}\left(  -1\right)  ^{i}\dbinom{n}{i}q^{i}%
\end{align*}
(here, we have renamed the summation index $k$ as $i$). Hence,%
\[
\sum_{i=0}^{n}\left(  -1\right)  ^{i}\dbinom{n}{i}q^{i}=\left(
\underbrace{\left(  -q\right)  +1}_{=1-q}\right)  ^{n}=\left(  1-q\right)
^{n}.
\]
This proves Corollary \ref{cor.binom.(1-q)n}.
\end{proof}

Let us next derive an easy corollary from Corollary \ref{cor.binom.1-1cor}:

\begin{corollary}
\label{cor.binom.1-1cor2}Let $n\in\mathbb{N}$. Let $i\in\mathbb{N}$. Then,%
\[
\sum_{j=0}^{n}\left(  -1\right)  ^{j}\dbinom{n}{j}\dbinom{j}{i}=\left(
-1\right)  ^{i}\left[  n=i\right]  .
\]

\end{corollary}

\begin{proof}
[Proof of Corollary \ref{cor.binom.1-1cor2}.]Each $j\in\left\{  0,1,\ldots
,i-1\right\}  $ satisfies%
\begin{equation}
\dbinom{j}{i}=0 \label{pf.cor.binom.1-1cor2.ji=0}%
\end{equation}
\footnote{\textit{Proof of (\ref{pf.cor.binom.1-1cor2.ji=0}):} Let
$j\in\left\{  0,1,\ldots,i-1\right\}  $. Thus, $j\leq i-1<i$ and $j\in\left\{
0,1,\ldots,i-1\right\}  \subseteq\mathbb{N}$. Hence, Proposition
\ref{prop.binom.0} (applied to $j$ and $i$ instead of $m$ and $n$) shows that
$\dbinom{j}{i}=0$. This proves (\ref{pf.cor.binom.1-1cor2.ji=0}).}.

We are in one of the following two cases:

\textit{Case 1:} We have $i\leq n$.

\textit{Case 2:} We have $i>n$.

Let us first consider Case 1. In this case, we have $i\leq n$. Thus,
$i\in\left\{  0,1,\ldots,n\right\}  $. Now,
\begin{align*}
&  \sum_{j=0}^{n}\left(  -1\right)  ^{j}\dbinom{n}{j}\dbinom{j}{i}\\
&  =\sum_{j=0}^{i-1}\left(  -1\right)  ^{j}\dbinom{n}{j}\underbrace{\dbinom
{j}{i}}_{\substack{=0\\\text{(by (\ref{pf.cor.binom.1-1cor2.ji=0}))}}%
}+\sum_{j=i}^{n}\underbrace{\left(  -1\right)  ^{j}}_{\substack{=\left(
-1\right)  ^{i+\left(  j+i\right)  }\\\text{(since }j\equiv2i+j=i+\left(
j+i\right)  \operatorname{mod}2\text{)}}}\dbinom{n}{j}\dbinom{j}{i}\\
&  \ \ \ \ \ \ \ \ \ \ \left(  \text{here, we have split the sum at
}j=i\text{, since }0\leq i\leq n\right) \\
&  =\underbrace{\sum_{j=0}^{i-1}\left(  -1\right)  ^{j}\dbinom{n}{j}0}%
_{=0}+\sum_{j=i}^{n}\left(  -1\right)  ^{i+\left(  j+i\right)  }\dbinom{n}%
{j}\dbinom{j}{i}=\sum_{j=i}^{n}\underbrace{\left(  -1\right)  ^{i+\left(
j+i\right)  }}_{=\left(  -1\right)  ^{i}\left(  -1\right)  ^{j+i}}\dbinom
{n}{j}\dbinom{j}{i}\\
&  =\sum_{j=i}^{n}\left(  -1\right)  ^{i}\left(  -1\right)  ^{j+i}\dbinom
{n}{j}\dbinom{j}{i}=\left(  -1\right)  ^{i}\underbrace{\sum_{j=i}^{n}\left(
-1\right)  ^{j+i}\dbinom{n}{j}\dbinom{j}{i}}_{\substack{=\left[  i=n\right]
\\\text{(by Corollary \ref{cor.binom.1-1cor})}}}\\
&  =\left(  -1\right)  ^{i}\left[  \underbrace{i=n}_{\Longleftrightarrow
\ \left(  n=i\right)  }\right]  =\left(  -1\right)  ^{i}\left[  n=i\right]  .
\end{align*}
Hence, Corollary \ref{cor.binom.1-1cor2} is proven in Case 1.

\begin{vershort}
Let us now consider Case 2. In this case, we have $i>n$. Thus, $n<i$, so that
$n\in\left\{  0,1,\ldots,i-1\right\}  $ (since $n\in\mathbb{N}$). But we don't
have $n=i$ (since we have $n<i$); thus, we have $\left[  n=i\right]  =0$.
Hence, $\left(  -1\right)  ^{i}\underbrace{\left[  n=i\right]  }_{=0}=0$. But
\[
\sum_{j=0}^{n}\left(  -1\right)  ^{j}\dbinom{n}{j}\underbrace{\dbinom{j}{i}%
}_{\substack{=0\\\text{(by (\ref{pf.cor.binom.1-1cor2.ji=0})}\\\text{(since
}j\in\left\{  0,1,\ldots,i-1\right\}  \\\text{(because }j\leq n<i\text{ and
}j\in\mathbb{N}\text{)))}}}=\sum_{j=0}^{n}\left(  -1\right)  ^{j}\dbinom{n}%
{j}0=0=\left(  -1\right)  ^{i}\left[  n=i\right]  .
\]
Hence, Corollary \ref{cor.binom.1-1cor2} is proven in Case 2.
\end{vershort}

\begin{verlong}
Let us now consider Case 2. In this case, we have $i>n$. Thus, $n<i$, so that
$n\leq i-1$ (since $n$ and $i$ are integers). Thus, $n\in\left\{
0,1,\ldots,i-1\right\}  $ (since $n\in\mathbb{N}$). But we don't have $n=i$
(since we have $n<i$); thus, we have $\left[  n=i\right]  =0$. Hence, $\left(
-1\right)  ^{i}\underbrace{\left[  n=i\right]  }_{=0}=0$. But each
$j\in\left\{  0,1,\ldots,n\right\}  $ satisfies%
\begin{equation}
\dbinom{j}{i}=0 \label{pf.cor.binom.1-1cor2.3}%
\end{equation}
\footnote{\textit{Proof of (\ref{pf.cor.binom.1-1cor2.3}):} Let $j\in\left\{
0,1,\ldots,n\right\}  $. Thus, $j\leq n<i$, so that $j\leq i-1$ (since $j$ and
$i$ are integers). Also, $j\in\left\{  0,1,\ldots,n\right\}  \subseteq
\mathbb{N}$. Hence, from $j\leq i-1$, we obtain $j\in\left\{  0,1,\ldots
,i-1\right\}  $. Thus, (\ref{pf.cor.binom.1-1cor2.ji=0}) shows that
$\dbinom{j}{i}=0$. This proves (\ref{pf.cor.binom.1-1cor2.3}).}. Now,%
\[
\sum_{j=0}^{n}\left(  -1\right)  ^{j}\dbinom{n}{j}\underbrace{\dbinom{j}{i}%
}_{\substack{=0\\\text{(by (\ref{pf.cor.binom.1-1cor2.3}))}}}=\sum_{j=0}%
^{n}\left(  -1\right)  ^{j}\dbinom{n}{j}0=0=\left(  -1\right)  ^{i}\left[
n=i\right]  .
\]
Hence, Corollary \ref{cor.binom.1-1cor2} is proven in Case 2.
\end{verlong}

We have now proven Corollary \ref{cor.binom.1-1cor2} in both Cases 1 and 2.
Hence, Corollary \ref{cor.binom.1-1cor2} always holds.
\end{proof}

Another simple lemma shall be of use to us:

\begin{lemma}
\label{lem.sol.bininv.pendulum}Let $n\in\mathbb{N}$.

\textbf{(a)} We have $\dfrac{1}{2}\left(  1+\left(  -1\right)  ^{n}\right)  =%
\begin{cases}
1, & \text{if }n\text{ is even;}\\
0, & \text{if }n\text{ is odd}%
\end{cases}
$.

\textbf{(b)} We have $\dfrac{1}{2}\left(  0^{n}+2^{n}\right)  =%
\begin{cases}
1, & \text{if }n=0;\\
2^{n-1}, & \text{if }n>0
\end{cases}
$.
\end{lemma}

\begin{vershort}
\begin{proof}
[Proof of Lemma \ref{lem.sol.bininv.pendulum}.]Each of the two claims of Lemma
\ref{lem.sol.bininv.pendulum} follows by a straightforward case distinction
(which we leave to the reader).
\end{proof}
\end{vershort}

\begin{verlong}
\begin{proof}
[Proof of Lemma \ref{lem.sol.bininv.pendulum}.]\textbf{(a)} We are in one of
the following two cases:

\textit{Case 1:} The integer $n$ is even.

\textit{Case 2:} The integer $n$ is odd.

Let us first consider Case 1. In this case, the integer $n$ is even. Thus,
$\left(  -1\right)  ^{n}=1$. Hence, $\dfrac{1}{2}\left(  1+\underbrace{\left(
-1\right)  ^{n}}_{=1}\right)  =\dfrac{1}{2}\left(  1+1\right)  =1$. Comparing
this with%
\[%
\begin{cases}
1, & \text{if }n\text{ is even;}\\
0, & \text{if }n\text{ is odd}%
\end{cases}
=1\ \ \ \ \ \ \ \ \ \ \left(  \text{since }n\text{ is even}\right)  ,
\]
we obtain $\dfrac{1}{2}\left(  1+\left(  -1\right)  ^{n}\right)  =%
\begin{cases}
1, & \text{if }n\text{ is even;}\\
0, & \text{if }n\text{ is odd}%
\end{cases}
$. Hence, Lemma \ref{lem.sol.bininv.pendulum} \textbf{(a)} holds in Case 1.

Let us now consider Case 2. In this case, the integer $n$ is odd. Thus,
$\left(  -1\right)  ^{n}=-1$. Hence, $\dfrac{1}{2}\left(
1+\underbrace{\left(  -1\right)  ^{n}}_{=-1}\right)  =\dfrac{1}{2}\left(
1+\left(  -1\right)  \right)  =0$. Comparing this with%
\[%
\begin{cases}
1, & \text{if }n\text{ is even;}\\
0, & \text{if }n\text{ is odd}%
\end{cases}
=0\ \ \ \ \ \ \ \ \ \ \left(  \text{since }n\text{ is odd}\right)  ,
\]
we obtain $\dfrac{1}{2}\left(  1+\left(  -1\right)  ^{n}\right)  =%
\begin{cases}
1, & \text{if }n\text{ is even;}\\
0, & \text{if }n\text{ is odd}%
\end{cases}
$. Hence, Lemma \ref{lem.sol.bininv.pendulum} \textbf{(a)} holds in Case 2.

We thus have proven Lemma \ref{lem.sol.bininv.pendulum} \textbf{(a)} in each
of the two Cases 1 and 2. Since these two Cases cover all possibilities, we
thus conclude that Lemma \ref{lem.sol.bininv.pendulum} \textbf{(a)} always holds.

\textbf{(b)} We have $n\geq0$ (since $n\in\mathbb{N}$). Hence, we are in one
of the following two cases:

\textit{Case 1:} We have $n=0$.

\textit{Case 2:} We have $n>0$.

Let us first consider Case 1. In this case, we have $n=0$. Thus, $0^{n}%
=0^{0}=1$. Also, from $n=0$, we obtain $2^{n}=2^{0}=1$. Thus, $\dfrac{1}%
{2}\left(  \underbrace{0^{n}}_{=1}+\underbrace{2^{n}}_{=1}\right)  =\dfrac
{1}{2}\left(  1+1\right)  =1$. Comparing this with%
\[%
\begin{cases}
1, & \text{if }n=0;\\
2^{n-1}, & \text{if }n>0
\end{cases}
=1\ \ \ \ \ \ \ \ \ \ \left(  \text{since }n=0\right)  ,
\]
we obtain $\dfrac{1}{2}\left(  0^{n}+2^{n}\right)  =%
\begin{cases}
1, & \text{if }n=0;\\
2^{n-1}, & \text{if }n>0
\end{cases}
$. Hence, Lemma \ref{lem.sol.bininv.pendulum} \textbf{(b)} holds in Case 1.

Let us now consider Case 2. In this case, we have $n>0$. Thus, $0^{n}=0$.
Hence, $\dfrac{1}{2}\left(  \underbrace{0^{n}}_{=0}+2^{n}\right)  =\dfrac
{1}{2}2^{n}=2^{n-1}$. Comparing this with%
\[%
\begin{cases}
1, & \text{if }n=0;\\
2^{n-1}, & \text{if }n>0
\end{cases}
=2^{n-1}\ \ \ \ \ \ \ \ \ \ \left(  \text{since }n>0\right)  ,
\]
we obtain $\dfrac{1}{2}\left(  0^{n}+2^{n}\right)  =%
\begin{cases}
1, & \text{if }n=0;\\
2^{n-1}, & \text{if }n>0
\end{cases}
$. Hence, Lemma \ref{lem.sol.bininv.pendulum} \textbf{(b)} holds in Case 2.

We thus have proven Lemma \ref{lem.sol.bininv.pendulum} \textbf{(b)} in each
of the two Cases 1 and 2. Since these two Cases cover all possibilities, we
thus conclude that Lemma \ref{lem.sol.bininv.pendulum} \textbf{(b)} always holds.
\end{proof}
\end{verlong}

We can now give answers to parts \textbf{(b)}, \textbf{(c)}, \textbf{(d)} and
\textbf{(e)} of Exercise \ref{exe.bininv}:

\begin{proposition}
\label{prop.sol.bininv.c}Let $a\in\mathbb{N}$. The binomial transform of the
list $\left(  \dbinom{0}{a},\dbinom{1}{a},\ldots,\dbinom{N}{a}\right)  $ is
the list
\[
\left(  \left(  -1\right)  ^{a}\left[  0=a\right]  ,\left(  -1\right)
^{a}\left[  1=a\right]  ,\ldots,\left(  -1\right)  ^{a}\left[  N=a\right]
\right)  .
\]
\footnotemark
\end{proposition}

\footnotetext{This list $\left(  \left(  -1\right)  ^{a}\left[  0=a\right]
,\left(  -1\right)  ^{a}\left[  1=a\right]  ,\ldots,\left(  -1\right)
^{a}\left[  N=a\right]  \right)  $ actually has a very simple form:
\par
\begin{itemize}
\item If it has an $\left(  a+1\right)  $-st entry (i.e., if $a\leq N$), then
this entry is $\left(  -1\right)  ^{a}\underbrace{\left[  a=a\right]
}_{\substack{=1\\\text{(since }a=a\text{)}}}=\left(  -1\right)  ^{a}$.
\par
\item All other entries are $0$ (because if $n\in\left\{  0,1,\ldots
,N\right\}  $ is such that $n\neq a$, then $\left(  -1\right)  ^{a}%
\underbrace{\left[  n=a\right]  }_{\substack{=0\\\text{(since }n\neq
a\text{)}}}=0$).
\end{itemize}
\par
Thus, this list $\left(  \left(  -1\right)  ^{a}\left[  0=a\right]  ,\left(
-1\right)  ^{a}\left[  1=a\right]  ,\ldots,\left(  -1\right)  ^{a}\left[
N=a\right]  \right)  $ can be rewritten as $\left(  0,0,\ldots,0,\left(
-1\right)  ^{a},0,0,\ldots,0\right)  $ (with the $\left(  -1\right)  ^{a}$
entry being placed in the $\left(  a+1\right)  $-st position) when $a\leq N$,
and as $\left(  0,0,\ldots,0\right)  $ if $a>N$.}

\begin{proof}
[Proof of Proposition \ref{prop.sol.bininv.c}.]Let $\left(  b_{0},b_{1}%
,\ldots,b_{N}\right)  $ be the binomial transform of the list \newline$\left(
\dbinom{0}{a},\dbinom{1}{a},\ldots,\dbinom{N}{a}\right)  $. Thus,%
\[
b_{n}=\sum_{i=0}^{n}\left(  -1\right)  ^{i}\dbinom{n}{i}\dbinom{i}%
{a}\ \ \ \ \ \ \ \ \ \ \text{for each }n\in\left\{  0,1,\ldots,N\right\}
\]
(by the definition of the binomial transform).

Hence, each $n\in\left\{  0,1,\ldots,N\right\}  $ satisfies
\begin{align*}
b_{n}  &  =\sum_{i=0}^{n}\left(  -1\right)  ^{i}\dbinom{n}{i}\dbinom{i}%
{a}=\sum_{j=0}^{n}\left(  -1\right)  ^{j}\dbinom{n}{j}\dbinom{j}{a}\\
&  \ \ \ \ \ \ \ \ \ \ \left(  \text{here, we have renamed the summation index
}i\text{ as }j\right) \\
&  =\left(  -1\right)  ^{a}\left[  n=a\right]  \ \ \ \ \ \ \ \ \ \ \left(
\text{by Corollary \ref{cor.binom.1-1cor2} (applied to }i=a\text{)}\right)  .
\end{align*}
In other words,%
\begin{equation}
\left(  b_{0},b_{1},\ldots,b_{N}\right)  =\left(  \left(  -1\right)
^{a}\left[  0=a\right]  ,\left(  -1\right)  ^{a}\left[  1=a\right]
,\ldots,\left(  -1\right)  ^{a}\left[  N=a\right]  \right)  .
\label{pf.prop.sol.bininv.c.1}%
\end{equation}

But recall that the binomial transform of the list $\left(  \dbinom{0}%
{a},\dbinom{1}{a},\ldots,\dbinom{N}{a}\right)  $ is \newline$\left(
b_{0},b_{1},\ldots,b_{N}\right)  $. In view of (\ref{pf.prop.sol.bininv.c.1}),
this rewrites as follows: The binomial transform of the list $\left(
\dbinom{0}{a},\dbinom{1}{a},\ldots,\dbinom{N}{a}\right)  $ is $\left(  \left(
-1\right)  ^{a}\left[  0=a\right]  ,\left(  -1\right)  ^{a}\left[  1=a\right]
,\ldots,\left(  -1\right)  ^{a}\left[  N=a\right]  \right)  $. This proves
Proposition \ref{prop.sol.bininv.c}.
\end{proof}

\begin{proposition}
\label{prop.sol.bininv.b}The binomial transform of the list $\left(
1,1,\ldots,1\right)  $ (with $N+1$ entries) is the list $\left(
1,0,0,\ldots,0\right)  $ (with one $1$ and $N$ zeroes).
\end{proposition}

\begin{proof}
[Proof of Proposition \ref{prop.sol.bininv.b}.]We have $\dbinom{n}{0}=1$ for
each $n\in\left\{  0,1,\ldots,N\right\}  $ (by (\ref{eq.binom.00}) (applied to
$m=n$)). Hence,
\begin{equation}
\left(  \dbinom{0}{0},\dbinom{1}{0},\ldots,\dbinom{N}{0}\right)  =\left(
1,1,\ldots,1\right)  \label{pf.prop.sol.bininv.b.1}%
\end{equation}
(with $N+1$ entries).

But Proposition \ref{prop.sol.bininv.c} (applied to $a=0$) shows that the
binomial transform of the list $\left(  \dbinom{0}{0},\dbinom{1}{0}%
,\ldots,\dbinom{N}{0}\right)  $ is the list
\begin{align*}
&  \left(  \underbrace{\left(  -1\right)  ^{0}}_{=1}\left[  0=0\right]
,\underbrace{\left(  -1\right)  ^{0}}_{=1}\left[  1=0\right]  ,\ldots
,\underbrace{\left(  -1\right)  ^{0}}_{=1}\left[  N=0\right]  \right) \\
&  =\left(  \left[  0=0\right]  ,\left[  1=0\right]  ,\ldots,\left[
N=0\right]  \right)  =\left(  1,0,0,\ldots,0\right)
\end{align*}
(with one $1$ and $N$ zeroes). In view of (\ref{pf.prop.sol.bininv.b.1}), this
rewrites as follows: The binomial transform of the list $\left(
1,1,\ldots,1\right)  $ (with $N+1$ entries) is the list $\left(
1,0,0,\ldots,0\right)  $ (with one $1$ and $N$ zeroes). This proves
Proposition \ref{prop.sol.bininv.b}.
\end{proof}

\begin{proposition}
\label{prop.sol.bininv.d}Let $q\in\mathbb{Z}$. The binomial transform of the
list $\left(  q^{0},q^{1},\ldots,q^{N}\right)  $ is $\left(  \left(
1-q\right)  ^{0},\left(  1-q\right)  ^{1},\ldots,\left(  1-q\right)
^{N}\right)  $.
\end{proposition}

\begin{proof}
[Proof of Proposition \ref{prop.sol.bininv.d}.]Let $\left(  b_{0},b_{1}%
,\ldots,b_{N}\right)  $ be the binomial transform of the list $\left(
q^{0},q^{1},\ldots,q^{N}\right)  $. Thus,%
\[
b_{n}=\sum_{i=0}^{n}\left(  -1\right)  ^{i}\dbinom{n}{i}q^{i}%
\ \ \ \ \ \ \ \ \ \ \text{for each }n\in\left\{  0,1,\ldots,N\right\}
\]
(by the definition of the binomial transform). Hence, each $n\in\left\{
0,1,\ldots,N\right\}  $ satisfies%
\[
b_{n}=\sum_{i=0}^{n}\left(  -1\right)  ^{i}\dbinom{n}{i}q^{i}=\left(
1-q\right)  ^{n}%
\]
(by Corollary \ref{cor.binom.(1-q)n}). In other words,
\begin{equation}
\left(  b_{0},b_{1},\ldots,b_{N}\right)  =\left(  \left(  1-q\right)
^{0},\left(  1-q\right)  ^{1},\ldots,\left(  1-q\right)  ^{N}\right)  .
\label{pf.prop.sol.bininv.d.1}%
\end{equation}

Recall that the binomial transform of the list $\left(  q^{0},q^{1}%
,\ldots,q^{N}\right)  $ is $\left(  b_{0},b_{1},\ldots,b_{N}\right)  $. In
view of (\ref{pf.prop.sol.bininv.d.1}), this rewrites as follows: The binomial
transform of the list $\left(  q^{0},q^{1},\ldots,q^{N}\right)  $ is $\left(
\left(  1-q\right)  ^{0},\left(  1-q\right)  ^{1},\ldots,\left(  1-q\right)
^{N}\right)  $. This proves Proposition \ref{prop.sol.bininv.d}.
\end{proof}

\begin{proposition}
\label{prop.sol.bininv.e}The binomial transform of the list $\left(
1,0,1,0,1,0,\ldots\right)  $ (with $N+1$ entries) is $\left(  1,2^{0}%
,2^{1},\ldots,2^{N-1}\right)  $.
\end{proposition}

\begin{proof}
[Proof of Proposition \ref{prop.sol.bininv.e}.]Let $\left(  a_{0},a_{1}%
,\ldots,a_{N}\right)  $ be the list $\left(  1,0,1,0,1,0,\ldots\right)  $
(with $N+1$ entries). Thus, for each $i\in\left\{  0,1,\ldots,N\right\}  $, we
have%
\begin{equation}
a_{i}=%
\begin{cases}
1, & \text{if }i\text{ is even;}\\
0, & \text{if }i\text{ is odd}%
\end{cases}
. \label{pf.prop.sol.bininv.e.1}%
\end{equation}

Each $i\in\mathbb{N}$ satisfies%
\begin{equation}
\dfrac{1}{2}\left(  1+\left(  -1\right)  ^{i}\right)  =%
\begin{cases}
1, & \text{if }i\text{ is even;}\\
0, & \text{if }i\text{ is odd}%
\end{cases}
\label{pf.prop.sol.bininv.e.2}%
\end{equation}
(by Lemma \ref{lem.sol.bininv.pendulum} \textbf{(a)} (applied to $n=i$)).
Hence, each $i\in\left\{  0,1,\ldots,N\right\}  $ satisfies%
\begin{align}
a_{i}  &  =%
\begin{cases}
1, & \text{if }i\text{ is even;}\\
0, & \text{if }i\text{ is odd}%
\end{cases}
\ \ \ \ \ \ \ \ \ \ \left(  \text{by (\ref{pf.prop.sol.bininv.e.1})}\right)
\nonumber\\
&  =\dfrac{1}{2}\left(  1+\left(  -1\right)  ^{i}\right)
\ \ \ \ \ \ \ \ \ \ \left(  \text{by (\ref{pf.prop.sol.bininv.e.2})}\right)  .
\label{pf.prop.sol.bininv.e.3}%
\end{align}

Let $\left(  b_{0},b_{1},\ldots,b_{N}\right)  $ be the binomial transform of
the list $\left(  1,0,1,0,1,0,\ldots\right)  $ (with $N+1$ entries). In other
words, $\left(  b_{0},b_{1},\ldots,b_{N}\right)  $ is the binomial transform
of the list $\left(  a_{0},a_{1},\ldots,a_{N}\right)  $ (because the list
$\left(  1,0,1,0,1,0,\ldots\right)  $ (with $N+1$ entries) is precisely
$\left(  a_{0},a_{1},\ldots,a_{N}\right)  $). Hence,%
\[
b_{n}=\sum_{i=0}^{n}\left(  -1\right)  ^{i}\dbinom{n}{i}a_{i}%
\ \ \ \ \ \ \ \ \ \ \text{for each }n\in\left\{  0,1,\ldots,N\right\}
\]
(by the definition of the binomial transform). Thus, for each $n\in\left\{
0,1,\ldots,N\right\}  $, we obtain%
\begin{align*}
b_{n}  &  =\sum_{i=0}^{n}\left(  -1\right)  ^{i}\dbinom{n}{i}\underbrace{a_{i}%
}_{\substack{=\dfrac{1}{2}\left(  1+\left(  -1\right)  ^{i}\right)
\\\text{(by (\ref{pf.prop.sol.bininv.e.3}))}}}=\sum_{i=0}^{n}%
\underbrace{\left(  -1\right)  ^{i}\dbinom{n}{i}\cdot\dfrac{1}{2}\left(
1+\left(  -1\right)  ^{i}\right)  }_{\substack{=\dfrac{1}{2}\left(  -1\right)
^{i}\dbinom{n}{i}\left(  1+\left(  -1\right)  ^{i}\right)  \\=\dfrac{1}%
{2}\left(  -1\right)  ^{i}\dbinom{n}{i}\cdot1+\dfrac{1}{2}\left(  -1\right)
^{i}\dbinom{n}{i}\left(  -1\right)  ^{i}}}\\
&  =\sum_{i=0}^{n}\left(  \dfrac{1}{2}\left(  -1\right)  ^{i}\dbinom{n}%
{i}\cdot1+\dfrac{1}{2}\left(  -1\right)  ^{i}\dbinom{n}{i}\left(  -1\right)
^{i}\right) \\
&  =\dfrac{1}{2}\sum_{i=0}^{n}\left(  -1\right)  ^{i}\dbinom{n}{i}%
\cdot\underbrace{1}_{=1^{i}}+\dfrac{1}{2}\sum_{i=0}^{n}\left(  -1\right)
^{i}\dbinom{n}{i}\left(  -1\right)  ^{i}\\
&  =\dfrac{1}{2}\underbrace{\sum_{i=0}^{n}\left(  -1\right)  ^{i}\dbinom{n}%
{i}\cdot1^{i}}_{\substack{=\left(  1-1\right)  ^{n}\\\text{(by Corollary
\ref{cor.binom.(1-q)n}, applied to }q=1\text{)}}}+\dfrac{1}{2}\underbrace{\sum
_{i=0}^{n}\left(  -1\right)  ^{i}\dbinom{n}{i}\left(  -1\right)  ^{i}%
}_{\substack{=\left(  1-\left(  -1\right)  \right)  ^{n}\\\text{(by Corollary
\ref{cor.binom.(1-q)n}, applied to }q=-1\text{)}}}\\
&  =\dfrac{1}{2}\left(  \underbrace{1-1}_{=0}\right)  ^{n}+\dfrac{1}{2}\left(
\underbrace{1-\left(  -1\right)  }_{=2}\right)  ^{n}=\dfrac{1}{2}0^{n}%
+\dfrac{1}{2}2^{n}\\
&  =\dfrac{1}{2}\left(  0^{n}+2^{n}\right)  =%
\begin{cases}
1, & \text{if }n=0;\\
2^{n-1}, & \text{if }n>0
\end{cases}
\end{align*}
(by Lemma \ref{lem.sol.bininv.pendulum} \textbf{(b)}). In other words,
\begin{equation}
\left(  b_{0},b_{1},\ldots,b_{N}\right)  =\left(  1,2^{0},2^{1},\ldots
,2^{N-1}\right)  . \label{pf.prop.sol.bininv.e.7}%
\end{equation}

Recall that the binomial transform of the list $\left(  1,0,1,0,1,0,\ldots
\right)  $ (with $N+1$ entries) is $\left(  b_{0},b_{1},\ldots,b_{N}\right)
$. In view of (\ref{pf.prop.sol.bininv.e.7}), this rewrites as follows: The
binomial transform of the list $\left(  1,0,1,0,1,0,\ldots\right)  $ (with
$N+1$ entries) is $\left(  1,2^{0},2^{1},\ldots,2^{N-1}\right)  $. This proves
Proposition \ref{prop.sol.bininv.e}.
\end{proof}

Next, we shall arm ourselves with another elementary lemma:

\begin{lemma}
\label{lem.sol.bininv.3.2}Let $N$, $n$ and $j$ be nonnegative integers such
that $N\geq n$ and $N\geq j$. Then,%
\[
\sum_{i=0}^{n}\left(  -1\right)  ^{i}\dbinom{n}{i}\dbinom{N-i}{j}=\dbinom
{N-n}{N-j}.
\]

\end{lemma}

There are two ways to prove Lemma \ref{lem.sol.bininv.3.2}. One way is
combinatorial (using the principle of inclusion and exclusion) and is
explained in \cite[proof of Identity 17.1]{Galvin}. The other way is
algebraic; this is the way we shall now present.

\begin{proof}
[Proof of Lemma \ref{lem.sol.bininv.3.2}.]We have $j\leq N$ (since $N\geq j$).
Thus, Proposition \ref{prop.vandermonde.consequences} \textbf{(e)} (applied to
$j$, $n$ and $N$ instead of $x$, $y$ and $n$) yields%
\begin{equation}
\dbinom{n-j-1}{N-j}=\sum_{k=0}^{N}\left(  -1\right)  ^{k-j}\dbinom{k}%
{j}\dbinom{n}{N-k}. \label{pf.lem.sol.bininv.3.2.1}%
\end{equation}
But $N-j\geq0$ (since $N\geq j$), so that $N-j\in\mathbb{N}$. Hence,
Proposition \ref{prop.binom.upper-neg} (applied to $N-n$ and $N-j$ instead of
$m$ and $n$) yields%
\begin{align}
\dbinom{N-n}{N-j}  &  =\left(  -1\right)  ^{N-j}\underbrace{\dbinom{\left(
N-j\right)  -\left(  N-n\right)  -1}{N-j}}_{\substack{=\dbinom{n-j-1}%
{N-j}\\\text{(since }\left(  N-j\right)  -\left(  N-n\right)  -1=n-j-1\text{)}%
}}\nonumber\\
&  =\left(  -1\right)  ^{N-j}\underbrace{\dbinom{n-j-1}{N-j}}_{\substack{=\sum
_{k=0}^{N}\left(  -1\right)  ^{k-j}\dbinom{k}{j}\dbinom{n}{N-k}\\\text{(by
(\ref{pf.lem.sol.bininv.3.2.1}))}}}=\left(  -1\right)  ^{N-j}\sum_{k=0}%
^{N}\left(  -1\right)  ^{k-j}\dbinom{k}{j}\dbinom{n}{N-k}\nonumber\\
&  =\sum_{k=0}^{N}\underbrace{\left(  -1\right)  ^{N-j}\left(  -1\right)
^{k-j}}_{\substack{=\left(  -1\right)  ^{\left(  N-j\right)  +\left(
k-j\right)  }=\left(  -1\right)  ^{N-k}\\\text{(since }\left(  N-j\right)
+\left(  k-j\right)  =N+k-2j\equiv N+k\equiv N-k\operatorname{mod}2\text{)}%
}}\dbinom{k}{j}\dbinom{n}{N-k}\nonumber\\
&  =\sum_{k=0}^{N}\left(  -1\right)  ^{N-k}\dbinom{k}{j}\dbinom{n}{N-k}.
\label{pf.lem.sol.bininv.3.2.3}%
\end{align}

If $i$ is an integer satisfying $i\geq n+1$, then%
\begin{equation}
\dbinom{n}{i}=0 \label{pf.lem.sol.bininv.3.2.triv0}%
\end{equation}
\footnote{\textit{Proof of (\ref{pf.lem.sol.bininv.3.2.triv0}):} Let $i$ be an
integer satisfying $i\geq n+1$. Thus, $i\geq n+1>n\geq0$, so that
$i\in\mathbb{N}$. Also, from $i>n$, we obtain $n<i$. Hence, Proposition
\ref{prop.binom.0} (applied to $n$ and $i$ instead of $m$ and $n$) yields
$\dbinom{n}{i}=0$. This proves (\ref{pf.lem.sol.bininv.3.2.triv0}).}.

But $0\leq n\leq N$ (since $N\geq n$). Hence, we can split the sum $\sum
_{i=0}^{N}\left(  -1\right)  ^{i}\dbinom{n}{i}\dbinom{N-i}{j}$ at $i=n$. We
thus find%
\begin{align*}
&  \sum_{i=0}^{N}\left(  -1\right)  ^{i}\dbinom{n}{i}\dbinom{N-i}{j}\\
&  =\sum_{i=0}^{n}\left(  -1\right)  ^{i}\dbinom{n}{i}\dbinom{N-i}{j}%
+\sum_{i=n+1}^{N}\left(  -1\right)  ^{i}\underbrace{\dbinom{n}{i}%
}_{\substack{=0\\\text{(by (\ref{pf.lem.sol.bininv.3.2.triv0}))}}}\dbinom
{N-i}{j}\\
&  =\sum_{i=0}^{n}\left(  -1\right)  ^{i}\dbinom{n}{i}\dbinom{N-i}%
{j}+\underbrace{\sum_{i=n+1}^{N}\left(  -1\right)  ^{i}0\dbinom{N-i}{j}}%
_{=0}\\
&  =\sum_{i=0}^{n}\left(  -1\right)  ^{i}\dbinom{n}{i}\dbinom{N-i}{j}.
\end{align*}
Hence,%
\begin{align*}
&  \sum_{i=0}^{n}\left(  -1\right)  ^{i}\dbinom{n}{i}\dbinom{N-i}{j}\\
&  =\sum_{i=0}^{N}\left(  -1\right)  ^{i}\dbinom{n}{i}\dbinom{N-i}{j}%
=\sum_{k=0}^{N}\left(  -1\right)  ^{N-k}\dbinom{n}{N-k}\underbrace{\dbinom
{N-\left(  N-k\right)  }{j}}_{\substack{=\dbinom{k}{j}\\\text{(since
}N-\left(  N-k\right)  =k\text{)}}}\\
&  \ \ \ \ \ \ \ \ \ \ \left(  \text{here, we have substituted }N-k\text{ for
}i\text{ in the sum}\right) \\
&  =\sum_{k=0}^{N}\left(  -1\right)  ^{N-k}\underbrace{\dbinom{n}{N-k}%
\dbinom{k}{j}}_{=\dbinom{k}{j}\dbinom{n}{N-k}}=\sum_{k=0}^{N}\left(
-1\right)  ^{N-k}\dbinom{k}{j}\dbinom{n}{N-k}=\dbinom{N-n}{N-j}%
\end{align*}
(by (\ref{pf.lem.sol.bininv.3.2.3})). This proves Lemma
\ref{lem.sol.bininv.3.2}.
\end{proof}

Let us now make two definitions that are slight variations on the definitions
of $B$ and $W$ made in Exercise \ref{exe.bininv}:

\begin{definition}
\label{def.sol.bininv.BW}\textbf{(a)} Let $B:\mathbb{Q}^{N+1}\rightarrow
\mathbb{Q}^{N+1}$ be the map which sends every list $\left(  f_{0}%
,f_{1},\ldots,f_{N}\right)  \in\mathbb{Q}^{N+1}$ to its binomial transform
$\left(  g_{0},g_{1},\ldots,g_{N}\right)  \in\mathbb{Q}^{N+1}$.

\textbf{(b)} Let $W:\mathbb{Q}^{N+1}\rightarrow\mathbb{Q}^{N+1}$ be the map
which sends every list $\left(  f_{0},f_{1},\ldots,f_{N}\right)  \in
\mathbb{Q}^{N+1}$ to the list $\left(  \left(  -1\right)  ^{N}f_{N},\left(
-1\right)  ^{N}f_{N-1},\ldots,\left(  -1\right)  ^{N}f_{0}\right)
\in\mathbb{Q}^{N+1}$.
\end{definition}

The maps $B$ and $W$ we have just defined are \textbf{almost} the same as the
maps $B$ and $W$ from Exercise \ref{exe.bininv} \textbf{(f)}. The only
difference is that the former maps are defined on $\mathbb{Q}^{N+1}$ (and have
codomains $\mathbb{Q}^{N+1}$ as well), whereas the latter are defined on
$\mathbb{Z}^{N+1}$ (and have codomains $\mathbb{Z}^{N+1}$). Thus, the latter
maps are restrictions of the former maps. Hence, in order to prove that the
latter maps satisfy the two equalities $B\circ W\circ B=W\circ B\circ W$ and
$\left(  B\circ W\right)  ^{3}=\operatorname*{id}$ (as demanded by Exercise
\ref{exe.bininv} \textbf{(f)}), it is perfectly sufficient to show that the
former maps satisfy these two equalities. In other words, it is perfectly
sufficient to prove the following theorem:

\begin{theorem}
\label{thm.sol.bininv.f}The two maps $B$ and $W$ introduced in Definition
\ref{def.sol.bininv.BW} satisfy $B\circ W\circ B=W\circ B\circ W$ and $\left(
B\circ W\right)  ^{3}=\operatorname*{id}$.
\end{theorem}

Before we prove Theorem \ref{thm.sol.bininv.f}, let us show three simpler facts:

\begin{lemma}
\label{lem.sol.bininv.f.Nn}Let $\mathbf{f}=\left(  f_{0},f_{1},\ldots
,f_{N}\right)  $ and $\mathbf{g}=\left(  g_{0},g_{1},\ldots,g_{N}\right)  $ be
two lists such that $\mathbf{g}=B\left(  \mathbf{f}\right)  $.

\textbf{(a)} Then,%
\[
g_{n}=\sum_{i=0}^{n}\left(  -1\right)  ^{i}\dbinom{n}{i}f_{i}%
\ \ \ \ \ \ \ \ \ \ \text{for each }n\in\left\{  0,1,\ldots,N\right\}  .
\]

\textbf{(b)} Also,%
\[
g_{n}=\sum_{i=0}^{N}\left(  -1\right)  ^{i}\dbinom{n}{i}f_{i}%
\ \ \ \ \ \ \ \ \ \ \text{for each }n\in\left\{  0,1,\ldots,N\right\}  .
\]

\end{lemma}

\begin{proof}
[Proof of Lemma \ref{lem.sol.bininv.f.Nn}.]\textbf{(a)} We know that the map
$B$ sends the list $\left(  f_{0},f_{1},\ldots,f_{N}\right)  $ to its binomial
transform (by the definition of $B$). In other words, $B\left(  \left(
f_{0},f_{1},\ldots,f_{N}\right)  \right)  $ is the binomial transform of
$\left(  f_{0},f_{1},\ldots,f_{N}\right)  $. In other words, $\left(
g_{0},g_{1},\ldots,g_{N}\right)  $ is the binomial transform of $\left(
f_{0},f_{1},\ldots,f_{N}\right)  $ (since $B\left(  \underbrace{\left(
f_{0},f_{1},\ldots,f_{N}\right)  }_{=\mathbf{f}}\right)  =B\left(
\mathbf{f}\right)  =\mathbf{g}=\left(  g_{0},g_{1},\ldots,g_{N}\right)  $). In
other words, we have%
\[
g_{n}=\sum_{i=0}^{n}\left(  -1\right)  ^{i}\dbinom{n}{i}f_{i}%
\ \ \ \ \ \ \ \ \ \ \text{for each }n\in\left\{  0,1,\ldots,N\right\}
\]
(by the definition of the binomial transform). This proves Lemma
\ref{lem.sol.bininv.f.Nn} \textbf{(a)}.

\textbf{(b)} Let $n\in\left\{  0,1,\ldots,N\right\}  $. Thus, $0\leq n\leq N$.

Lemma \ref{lem.sol.bininv.f.Nn} \textbf{(a)} yields%
\begin{equation}
g_{n}=\sum_{i=0}^{n}\left(  -1\right)  ^{i}\dbinom{n}{i}f_{i}.
\label{pf.lem.sol.bininv.f.Nn.0}%
\end{equation}

For each integer $i$ satisfying $i\geq n+1$, we have%
\begin{equation}
\dbinom{n}{i}=0 \label{pf.lem.sol.bininv.f.Nn.1}%
\end{equation}
\footnote{\textit{Proof of (\ref{pf.lem.sol.bininv.f.Nn.1}):} Let $i$ be an
integer satisfying $i\geq n+1$. We have $i\geq n+1>n\geq0$, so that
$i\in\mathbb{N}$. Also, $n\in\left\{  0,1,\ldots,N\right\}  \subseteq
\mathbb{N}$. Also, $n<i$ (since $i>n$). Hence, Proposition \ref{prop.binom.0}
(applied to $n$ and $i$ instead of $m$ and $n$) yields $\dbinom{n}{i}=0$. This
proves (\ref{pf.lem.sol.bininv.f.Nn.1}).}.

But $0\leq n\leq N$. Hence, we can split the sum $\sum_{i=0}^{N}\left(
-1\right)  ^{i}\dbinom{n}{i}f_{i}$ at $i=n$. We thus find
\begin{align*}
\sum_{i=0}^{N}\left(  -1\right)  ^{i}\dbinom{n}{i}f_{i}  &  =\sum_{i=0}%
^{n}\left(  -1\right)  ^{i}\dbinom{n}{i}f_{i}+\sum_{i=n+1}^{N}\left(
-1\right)  ^{i}\underbrace{\dbinom{n}{i}}_{\substack{=0\\\text{(by
(\ref{pf.lem.sol.bininv.f.Nn.1}))}}}f_{i}\\
&  =\sum_{i=0}^{n}\left(  -1\right)  ^{i}\dbinom{n}{i}f_{i}+\underbrace{\sum
_{i=n+1}^{N}\left(  -1\right)  ^{i}0f_{i}}_{=0}=\sum_{i=0}^{n}\left(
-1\right)  ^{i}\dbinom{n}{i}f_{i}=g_{n}%
\end{align*}
(by (\ref{pf.lem.sol.bininv.f.Nn.0})). In other words, $g_{n}=\sum_{i=0}%
^{N}\left(  -1\right)  ^{i}\dbinom{n}{i}f_{i}$. This proves Lemma
\ref{lem.sol.bininv.f.Nn} \textbf{(b)}.
\end{proof}

\begin{proposition}
\label{prop.sol.bininv.f.BB}The map $B$ introduced in Definition
\ref{def.sol.bininv.BW} satisfies $B\circ B=\operatorname*{id}$.
\end{proposition}

\begin{proof}
[Proof of Proposition \ref{prop.sol.bininv.f.BB}.]Let $\mathbf{a}$ be a list.
Write this list $\mathbf{a}$ in the form $\mathbf{a}=\left(  a_{0}%
,a_{1},\ldots,a_{N}\right)  $.

Write the list $B\left(  \mathbf{a}\right)  $ in the form $B\left(
\mathbf{a}\right)  =\left(  b_{0},b_{1},\ldots,b_{N}\right)  $.

We know that the map $B$ sends the list $\left(  a_{0},a_{1},\ldots
,a_{N}\right)  $ to its binomial transform (by the definition of $B$). In
other words, $B\left(  \left(  a_{0},a_{1},\ldots,a_{N}\right)  \right)  $ is
the binomial transform of $\left(  a_{0},a_{1},\ldots,a_{N}\right)  $. In
other words, $\left(  b_{0},b_{1},\ldots,b_{N}\right)  $ is the binomial
transform of $\left(  a_{0},a_{1},\ldots,a_{N}\right)  $ (since $B\left(
\underbrace{\left(  a_{0},a_{1},\ldots,a_{N}\right)  }_{=\mathbf{a}}\right)
=B\left(  \mathbf{a}\right)  =\left(  b_{0},b_{1},\ldots,b_{N}\right)  $).
Hence, Corollary \ref{cor.sol.bininv.a} (applied to $f_{i}=a_{i}$ and
$g_{i}=b_{i}$) shows that $\left(  a_{0},a_{1},\ldots,a_{N}\right)  $ is the
binomial transform of $\left(  b_{0},b_{1},\ldots,b_{N}\right)  $.

We know that the map $B$ sends the list $\left(  b_{0},b_{1},\ldots
,b_{N}\right)  $ to its binomial transform (by the definition of $B$). In
other words, $B\left(  \left(  b_{0},b_{1},\ldots,b_{N}\right)  \right)  $ is
the binomial transform of $\left(  b_{0},b_{1},\ldots,b_{N}\right)  $.

The two lists $B\left(  \left(  b_{0},b_{1},\ldots,b_{N}\right)  \right)  $
and $\left(  a_{0},a_{1},\ldots,a_{N}\right)  $ must be equal, since each of
them is the binomial transform of $\left(  b_{0},b_{1},\ldots,b_{N}\right)  $
(as we have proven above). In other words, $B\left(  \left(  b_{0}%
,b_{1},\ldots,b_{N}\right)  \right)  =\left(  a_{0},a_{1},\ldots,a_{N}\right)
$. Now,%
\begin{align*}
\left(  B\circ B\right)  \left(  \mathbf{a}\right)   &  =B\left(
\underbrace{B\left(  \mathbf{a}\right)  }_{=\left(  b_{0},b_{1},\ldots
,b_{N}\right)  }\right)  =B\left(  \left(  b_{0},b_{1},\ldots,b_{N}\right)
\right)  =\left(  a_{0},a_{1},\ldots,a_{N}\right) \\
&  =\mathbf{a}=\operatorname*{id}\left(  \mathbf{a}\right)  .
\end{align*}

Now, forget that we fixed $\mathbf{a}$. We thus have shown that $\left(
B\circ B\right)  \left(  \mathbf{a}\right)  =\operatorname*{id}\left(
\mathbf{a}\right)  $ for each list $\mathbf{a}$. In other words, $B\circ
B=\operatorname*{id}$. This proves Proposition \ref{prop.sol.bininv.f.BB}.
\end{proof}

\begin{proposition}
\label{prop.sol.bininv.f.WW}The map $W$ introduced in Definition
\ref{def.sol.bininv.BW} satisfies $W\circ W=\operatorname*{id}$.
\end{proposition}

\begin{proof}
[Proof of Proposition \ref{prop.sol.bininv.f.WW}.]Let $\mathbf{a}$ be a list.
Write $\mathbf{a}$ in the form $\mathbf{a}=\left(  a_{0},a_{1},\ldots
,a_{N}\right)  $.

The map $W$ sends the list $\left(  a_{0},a_{1},\ldots,a_{N}\right)  $ to
$\left(  \left(  -1\right)  ^{N}a_{N},\left(  -1\right)  ^{N}a_{N-1}%
,\ldots,\left(  -1\right)  ^{N}a_{0}\right)  $ (by the definition of $W$). In
other words,%
\[
W\left(  \left(  a_{0},a_{1},\ldots,a_{N}\right)  \right)  =\left(  \left(
-1\right)  ^{N}a_{N},\left(  -1\right)  ^{N}a_{N-1},\ldots,\left(  -1\right)
^{N}a_{0}\right)  .
\]

Let $\mathbf{b}$ be the list $W\left(  \mathbf{a}\right)  $. Thus,
$\mathbf{b}=W\left(  \mathbf{a}\right)  $. Write the list $\mathbf{b}$ in the
form $\mathbf{b}=\left(  b_{0},b_{1},\ldots,b_{N}\right)  $. Thus,
\begin{align*}
\left(  b_{0},b_{1},\ldots,b_{N}\right)   &  =\mathbf{b}=W\left(
\underbrace{\mathbf{a}}_{=\left(  a_{0},a_{1},\ldots,a_{N}\right)  }\right) \\
&  =W\left(  \left(  a_{0},a_{1},\ldots,a_{N}\right)  \right)  =\left(
\left(  -1\right)  ^{N}a_{N},\left(  -1\right)  ^{N}a_{N-1},\ldots,\left(
-1\right)  ^{N}a_{0}\right)  .
\end{align*}
In other words,%
\begin{equation}
b_{n}=\left(  -1\right)  ^{N}a_{N-n}\ \ \ \ \ \ \ \ \ \ \text{for each }%
n\in\left\{  0,1,\ldots,N\right\}  . \label{pf.prop.sol.bininv.f.WW.1}%
\end{equation}

Hence, for each $n\in\left\{  0,1,\ldots,N\right\}  $, we have
\begin{align*}
\left(  -1\right)  ^{N}\underbrace{b_{N-n}}_{\substack{=\left(  -1\right)
^{N}a_{N-\left(  N-n\right)  }\\\text{(by (\ref{pf.prop.sol.bininv.f.WW.1})
(applied}\\\text{to }N-n\text{ instead of }n\text{))}}}  &
=\underbrace{\left(  -1\right)  ^{N}\left(  -1\right)  ^{N}}_{=\left(  \left(
-1\right)  \left(  -1\right)  \right)  ^{N}}\underbrace{a_{N-\left(
N-n\right)  }}_{\substack{=a_{n}\\\text{(since }N-\left(  N-n\right)
=n\text{)}}}\\
&  =\left(  \underbrace{\left(  -1\right)  \left(  -1\right)  }_{=1}\right)
^{N}a_{n}=\underbrace{1^{N}}_{=1}a_{n}=a_{n}.
\end{align*}
In other words,
\[
\left(  \left(  -1\right)  ^{N}b_{N},\left(  -1\right)  ^{N}b_{N-1}%
,\ldots,\left(  -1\right)  ^{N}b_{0}\right)  =\left(  a_{0},a_{1},\ldots
,a_{N}\right)  .
\]

But the map $W$ sends the list $\left(  b_{0},b_{1},\ldots,b_{N}\right)  $ to
$\left(  \left(  -1\right)  ^{N}b_{N},\left(  -1\right)  ^{N}b_{N-1}%
,\ldots,\left(  -1\right)  ^{N}b_{0}\right)  $ (by the definition of $W$). In
other words,%
\[
W\left(  \left(  b_{0},b_{1},\ldots,b_{N}\right)  \right)  =\left(  \left(
-1\right)  ^{N}b_{N},\left(  -1\right)  ^{N}b_{N-1},\ldots,\left(  -1\right)
^{N}b_{0}\right)  .
\]
Hence,%
\begin{align*}
\left(  W\circ W\right)  \left(  \mathbf{a}\right)   &  =W\left(
\underbrace{W\left(  \mathbf{a}\right)  }_{=\mathbf{b}=\left(  b_{0}%
,b_{1},\ldots,b_{N}\right)  }\right)  =W\left(  \left(  b_{0},b_{1}%
,\ldots,b_{N}\right)  \right) \\
&  =\left(  \left(  -1\right)  ^{N}b_{N},\left(  -1\right)  ^{N}b_{N-1}%
,\ldots,\left(  -1\right)  ^{N}b_{0}\right)  =\left(  a_{0},a_{1},\ldots
,a_{N}\right) \\
&  =\mathbf{a}=\operatorname*{id}\left(  \mathbf{a}\right)  .
\end{align*}

Now, forget that we fixed $\mathbf{a}$. We thus have shown that $\left(
W\circ W\right)  \left(  \mathbf{a}\right)  =\operatorname*{id}\left(
\mathbf{a}\right)  $ for each list $\mathbf{a}$. In other words, $W\circ
W=\operatorname*{id}$. This proves Proposition \ref{prop.sol.bininv.f.WW}.
\end{proof}

We are now ready to prove Theorem \ref{thm.sol.bininv.f}:

\begin{proof}
[Proof of Theorem \ref{thm.sol.bininv.f}.]Let us first focus on proving that
$B\circ W\circ B=W\circ B\circ W$.

Indeed, let $\mathbf{a}$ be a list. Write $\mathbf{a}$ in the form
$\mathbf{a}=\left(  a_{0},a_{1},\ldots,a_{N}\right)  $.

Let $\mathbf{b}$ be the list $B\left(  \mathbf{a}\right)  $. Thus,
$\mathbf{b}=B\left(  \mathbf{a}\right)  $. Write the list $\mathbf{b}$ in the
form $\mathbf{b}=\left(  b_{0},b_{1},\ldots,b_{N}\right)  $.

Lemma \ref{lem.sol.bininv.f.Nn} \textbf{(b)} (applied to $\mathbf{a}$, $a_{i}%
$, $\mathbf{b}$ and $b_{i}$ instead of $\mathbf{f}$, $f_{i}$, $\mathbf{g}$ and
$g_{i}$) shows that
\begin{equation}
b_{n}=\sum_{i=0}^{N}\left(  -1\right)  ^{i}\dbinom{n}{i}a_{i}%
\ \ \ \ \ \ \ \ \ \ \text{for each }n\in\left\{  0,1,\ldots,N\right\}
\label{sol.bininv.3.sol1.bn=n}%
\end{equation}
(since $\mathbf{b}=B\left(  \mathbf{a}\right)  $).

The map $W$ sends the list $\left(  b_{0},b_{1},\ldots,b_{N}\right)  $ to
$\left(  \left(  -1\right)  ^{N}b_{N},\left(  -1\right)  ^{N}b_{N-1}%
,\ldots,\left(  -1\right)  ^{N}b_{0}\right)  $ (by the definition of $W$). In
other words,%
\[
W\left(  \left(  b_{0},b_{1},\ldots,b_{N}\right)  \right)  =\left(  \left(
-1\right)  ^{N}b_{N},\left(  -1\right)  ^{N}b_{N-1},\ldots,\left(  -1\right)
^{N}b_{0}\right)  .
\]
Thus,%
\begin{align}
W\left(  \underbrace{\mathbf{b}}_{=\left(  b_{0},b_{1},\ldots,b_{N}\right)
}\right)   &  =W\left(  \left(  b_{0},b_{1},\ldots,b_{N}\right)  \right)
\nonumber\\
&  =\left(  \left(  -1\right)  ^{N}b_{N},\left(  -1\right)  ^{N}b_{N-1}%
,\ldots,\left(  -1\right)  ^{N}b_{0}\right)  . \label{sol.bininv.3.sol1.Wb=}%
\end{align}

Now, define a list $\mathbf{c}$ by $\mathbf{c}=B\left(  W\left(
\mathbf{b}\right)  \right)  $. Write the list $\mathbf{c}$ in the form
$\mathbf{c}=\left(  c_{0},c_{1},\ldots,c_{N}\right)  $. Hence, Lemma
\ref{lem.sol.bininv.f.Nn} \textbf{(a)} (applied to $W\left(  \mathbf{b}%
\right)  $, $\left(  -1\right)  ^{N}b_{N-i}$, $\mathbf{c}$ and $c_{i}$ instead
of $\mathbf{f}$, $f_{i}$, $\mathbf{g}$ and $g_{i}$) shows that
\begin{equation}
c_{n}=\sum_{i=0}^{n}\left(  -1\right)  ^{i}\dbinom{n}{i}\left(  -1\right)
^{N}b_{N-i}\ \ \ \ \ \ \ \ \ \ \text{for each }n\in\left\{  0,1,\ldots
,N\right\}  \label{sol.bininv.3.sol1.cn=}%
\end{equation}
(because $\mathbf{c}=B\left(  W\left(  \mathbf{b}\right)  \right)  $ and
because of (\ref{sol.bininv.3.sol1.Wb=})).

The map $W$ sends the list $\left(  a_{0},a_{1},\ldots,a_{N}\right)  $ to
$\left(  \left(  -1\right)  ^{N}a_{N},\left(  -1\right)  ^{N}a_{N-1}%
,\ldots,\left(  -1\right)  ^{N}a_{0}\right)  $ (by the definition of $W$). In
other words,%
\[
W\left(  \left(  a_{0},a_{1},\ldots,a_{N}\right)  \right)  =\left(  \left(
-1\right)  ^{N}a_{N},\left(  -1\right)  ^{N}a_{N-1},\ldots,\left(  -1\right)
^{N}a_{0}\right)  .
\]
Thus,%
\begin{align}
W\left(  \underbrace{\mathbf{a}}_{=\left(  a_{0},a_{1},\ldots,a_{N}\right)
}\right)   &  =W\left(  \left(  a_{0},a_{1},\ldots,a_{N}\right)  \right)
\nonumber\\
&  =\left(  \left(  -1\right)  ^{N}a_{N},\left(  -1\right)  ^{N}a_{N-1}%
,\ldots,\left(  -1\right)  ^{N}a_{0}\right)  . \label{sol.bininv.3.sol1.Wa=}%
\end{align}

Define a list $\mathbf{d}$ by $\mathbf{d}=B\left(  W\left(  \mathbf{a}\right)
\right)  $. Write the list $\mathbf{d}$ in the form $\mathbf{d}=\left(
d_{0},d_{1},\ldots,d_{N}\right)  $. Hence, Lemma \ref{lem.sol.bininv.f.Nn}
\textbf{(b)} (applied to $W\left(  \mathbf{a}\right)  $, $\left(  -1\right)
^{N}a_{N-i}$, $\mathbf{d}$ and $d_{i}$ instead of $\mathbf{f}$, $f_{i}$,
$\mathbf{g}$ and $g_{i}$) shows that%
\begin{equation}
d_{n}=\sum_{i=0}^{N}\left(  -1\right)  ^{i}\dbinom{n}{i}\left(  -1\right)
^{N}a_{N-i}\ \ \ \ \ \ \ \ \ \ \text{for each }n\in\left\{  0,1,\ldots
,N\right\}  \label{sol.bininv.3.sol1.dn=}%
\end{equation}
(because $\mathbf{d}=B\left(  W\left(  \mathbf{a}\right)  \right)  $ and
because of (\ref{sol.bininv.3.sol1.Wa=})).

The map $W$ sends the list $\left(  d_{0},d_{1},\ldots,d_{N}\right)  $ to
$\left(  \left(  -1\right)  ^{N}d_{N},\left(  -1\right)  ^{N}d_{N-1}%
,\ldots,\left(  -1\right)  ^{N}d_{0}\right)  $ (by the definition of $W$). In
other words,%
\[
W\left(  \left(  d_{0},d_{1},\ldots,d_{N}\right)  \right)  =\left(  \left(
-1\right)  ^{N}d_{N},\left(  -1\right)  ^{N}d_{N-1},\ldots,\left(  -1\right)
^{N}d_{0}\right)  .
\]
Thus,%
\begin{align}
W\left(  \underbrace{\mathbf{d}}_{=\left(  d_{0},d_{1},\ldots,d_{N}\right)
}\right)   &  =W\left(  \left(  d_{0},d_{1},\ldots,d_{N}\right)  \right)
\nonumber\\
&  =\left(  \left(  -1\right)  ^{N}d_{N},\left(  -1\right)  ^{N}d_{N-1}%
,\ldots,\left(  -1\right)  ^{N}d_{0}\right)  . \label{sol.bininv.3.sol1.Wd=}%
\end{align}

We shall now show that $\mathbf{c}=W\left(  \mathbf{d}\right)  $.

Indeed, for any $g\in\left\{  0,1,\ldots,N\right\}  $, we have%
\begin{align}
b_{g}  &  =\sum_{i=0}^{N}\left(  -1\right)  ^{i}\dbinom{g}{i}a_{i}%
\ \ \ \ \ \ \ \ \ \ \left(  \text{by (\ref{sol.bininv.3.sol1.bn=n}), applied
to }n=g\right) \nonumber\\
&  =\sum_{j=0}^{N}\left(  -1\right)  ^{j}\dbinom{g}{j}a_{j}
\label{pf.thm.sol.bininv.f.bg=}%
\end{align}
(here, we have renamed the summation index $i$ as $j$).

Now, let $n\in\left\{  0,1,\ldots,N\right\}  $ be arbitrary. Then, $n\leq N$,
so that $N\geq n$. Hence, $N-n\geq0$, so that $0\leq N-n\leq N$ (since
$n\geq0$ (since $n\in\left\{  0,1,\ldots,N\right\}  $)). Thus, $N-n\in\left\{
0,1,\ldots,N\right\}  $. Hence, (\ref{sol.bininv.3.sol1.dn=}) (applied to
$N-n$ instead of $n$) yields%
\[
d_{N-n}=\sum_{i=0}^{N}\left(  -1\right)  ^{i}\dbinom{N-n}{i}\left(  -1\right)
^{N}a_{N-i}.
\]
Multiplying both sides of this equality by $\left(  -1\right)  ^{N}$, we
obtain%
\begin{align}
\left(  -1\right)  ^{N}d_{N-n}  &  =\left(  -1\right)  ^{N}\left(  \sum
_{i=0}^{N}\left(  -1\right)  ^{i}\dbinom{N-n}{i}\left(  -1\right)  ^{N}%
a_{N-i}\right) \nonumber\\
&  =\sum_{i=0}^{N}\left(  -1\right)  ^{i}\dbinom{N-n}{i}\underbrace{\left(
-1\right)  ^{N}\left(  -1\right)  ^{N}}_{\substack{=\left(  \left(  -1\right)
\left(  -1\right)  \right)  ^{N}=1^{N}\\\text{(since }\left(  -1\right)
\left(  -1\right)  =1\text{)}}}a_{N-i}\nonumber\\
&  =\sum_{i=0}^{N}\left(  -1\right)  ^{i}\dbinom{N-n}{i}\underbrace{1^{N}%
}_{=1}a_{N-i}\nonumber\\
&  =\sum_{i=0}^{N}\left(  -1\right)  ^{i}\dbinom{N-n}{i}a_{N-i}.
\label{pf.thm.sol.bininv.f-1NdN-n=}%
\end{align}

For each $i\in\left\{  0,1,\ldots,n\right\}  $, we have%
\begin{equation}
b_{N-i}=\sum_{j=0}^{N}\left(  -1\right)  ^{j}\dbinom{N-i}{j}a_{j}
\label{pf.thm.sol.bininv.bN-i=}%
\end{equation}
\footnote{\textit{Proof of (\ref{pf.thm.sol.bininv.bN-i=}):} Let $i\in\left\{
0,1,\ldots,n\right\}  $. Thus, $i\in\left\{  0,1,\ldots,n\right\}
\subseteq\left\{  0,1,\ldots,N\right\}  $ (since $n\leq N$). Hence,
$N-i\in\left\{  0,1,\ldots,N\right\}  $. Thus, (\ref{pf.thm.sol.bininv.f.bg=})
(applied to $g=N-i$) yields $b_{N-i}=\sum_{j=0}^{N}\left(  -1\right)
^{j}\dbinom{N-i}{j}a_{j}$. This proves (\ref{pf.thm.sol.bininv.bN-i=}).}.

But (\ref{sol.bininv.3.sol1.cn=}) becomes%
\begin{align*}
c_{n}  &  =\sum_{i=0}^{n}\left(  -1\right)  ^{i}\dbinom{n}{i}\left(
-1\right)  ^{N}\underbrace{b_{N-i}}_{\substack{=\sum_{j=0}^{N}\left(
-1\right)  ^{j}\dbinom{N-i}{j}a_{j}\\\text{(by (\ref{pf.thm.sol.bininv.bN-i=}%
))}}}\\
&  =\sum_{i=0}^{n}\left(  -1\right)  ^{i}\dbinom{n}{i}\left(  -1\right)
^{N}\sum_{j=0}^{N}\left(  -1\right)  ^{j}\dbinom{N-i}{j}a_{j}\\
&  =\underbrace{\sum_{i=0}^{n}\sum_{j=0}^{N}}_{=\sum_{j=0}^{N}\sum_{i=0}^{n}%
}\left(  -1\right)  ^{i}\dbinom{n}{i}\left(  -1\right)  ^{N}\left(  -1\right)
^{j}\dbinom{N-i}{j}a_{j}\\
&  =\sum_{j=0}^{N}\sum_{i=0}^{n}\left(  -1\right)  ^{i}\dbinom{n}{i}\left(
-1\right)  ^{N}\left(  -1\right)  ^{j}\dbinom{N-i}{j}a_{j}\\
&  =\sum_{j=0}^{N}\underbrace{\left(  -1\right)  ^{N}\left(  -1\right)  ^{j}%
}_{\substack{=\left(  -1\right)  ^{N+j}=\left(  -1\right)  ^{N-j}%
\\\text{(since }N+j\equiv N-j\operatorname{mod}2\text{)}}}\underbrace{\left(
\sum_{i=0}^{n}\left(  -1\right)  ^{i}\dbinom{n}{i}\dbinom{N-i}{j}\right)
}_{\substack{=\dbinom{N-n}{N-j}\\\text{(by Lemma \ref{lem.sol.bininv.3.2}
(since }N\geq j\\\text{(since }j\leq N\text{)))}}}\underbrace{a_{j}%
}_{\substack{=a_{N-\left(  N-j\right)  }\\\text{(since }j=N-\left(
N-j\right)  \text{)}}}\\
&  =\sum_{j=0}^{N}\left(  -1\right)  ^{N-j}\dbinom{N-n}{N-j}a_{N-\left(
N-j\right)  }=\sum_{i=0}^{N}\left(  -1\right)  ^{i}\dbinom{N-n}{i}a_{N-i}\\
&  \ \ \ \ \ \ \ \ \ \ \left(  \text{here, we have substituted }i\text{ for
}N-j\text{ in the sum}\right) \\
&  =\left(  -1\right)  ^{N}d_{N-n}\ \ \ \ \ \ \ \ \ \ \left(  \text{by
(\ref{pf.thm.sol.bininv.f-1NdN-n=})}\right)  .
\end{align*}

Now, forget that we fixed $n$. We thus have proven that $c_{n}=\left(
-1\right)  ^{N}d_{N-n}$ for each $n\in\left\{  0,1,\ldots,N\right\}  $. In
other words,%
\[
\left(  c_{0},c_{1},\ldots,c_{N}\right)  =\left(  \left(  -1\right)  ^{N}%
d_{N},\left(  -1\right)  ^{N}d_{N-1},\ldots,\left(  -1\right)  ^{N}%
d_{0}\right)  .
\]
Thus,%
\begin{align*}
\left(  B\circ W\circ B\right)  \left(  \mathbf{a}\right)   &  =B\left(
W\left(  \underbrace{B\left(  \mathbf{a}\right)  }_{=\mathbf{b}}\right)
\right)  =B\left(  W\left(  \mathbf{b}\right)  \right)  =\mathbf{c}%
\ \ \ \ \ \ \ \ \ \ \left(  \text{since }\mathbf{c}=B\left(  W\left(
\mathbf{b}\right)  \right)  \right) \\
&  =\left(  c_{0},c_{1},\ldots,c_{N}\right)  =\left(  \left(  -1\right)
^{N}d_{N},\left(  -1\right)  ^{N}d_{N-1},\ldots,\left(  -1\right)  ^{N}%
d_{0}\right) \\
&  =W\left(  \underbrace{\mathbf{d}}_{=B\left(  W\left(  \mathbf{a}\right)
\right)  }\right)  \ \ \ \ \ \ \ \ \ \ \left(  \text{by
(\ref{sol.bininv.3.sol1.Wd=})}\right) \\
&  =W\left(  B\left(  W\left(  \mathbf{a}\right)  \right)  \right)  =\left(
W\circ B\circ W\right)  \left(  \mathbf{a}\right)  .
\end{align*}

Now, forget that we fixed $\mathbf{a}$. We thus have proven that
\newline$\left(  B\circ W\circ B\right)  \left(  \mathbf{a}\right)  =\left(
W\circ B\circ W\right)  \left(  \mathbf{a}\right)  $ for each list
$\mathbf{a}$. In other words,%
\begin{equation}
B\circ W\circ B=W\circ B\circ W. \label{sol.bininv.3.sol1.braid}%
\end{equation}

Hence,%
\begin{align*}
\left(  B\circ W\right)  ^{3}  &  =\left(  B\circ W\right)  \circ\left(
B\circ W\right)  \circ\left(  B\circ W\right) \\
&  =\underbrace{B\circ W\circ B}_{\substack{=W\circ B\circ W\\\text{(by
(\ref{sol.bininv.3.sol1.braid})) }}}\circ W\circ B\circ W=W\circ
B\circ\underbrace{W\circ W}_{\substack{=\operatorname*{id}\\\text{(by
Proposition \ref{prop.sol.bininv.f.WW})}}}\circ B\circ W\\
&  =W\circ\underbrace{B\circ B}_{\substack{=\operatorname*{id}\\\text{(by
Proposition \ref{prop.sol.bininv.f.BB})}}}\circ W=W\circ W=\operatorname*{id}%
\ \ \ \ \ \ \ \ \ \ \left(  \text{by Proposition \ref{prop.sol.bininv.f.WW}%
}\right)  .
\end{align*}
This completes the proof of Theorem \ref{thm.sol.bininv.f}.
\end{proof}

\begin{proof}
[Solution to Exercise \ref{exe.bininv}.]Recall that the definition of the
binomial transform in Definition \ref{def.sol.bininv.bintr} generalizes the
definition of the binomial transform we gave in Exercise \ref{exe.bininv}.
Hence, part \textbf{(a)} of Exercise \ref{exe.bininv} follows from Corollary
\ref{cor.sol.bininv.a}. Part \textbf{(b)} follows from Proposition
\ref{prop.sol.bininv.b}. Part \textbf{(c)} follows from Proposition
\ref{prop.sol.bininv.c}. Part \textbf{(d)} follows from Proposition
\ref{prop.sol.bininv.d}. Part \textbf{(e)} follows from Proposition
\ref{prop.sol.bininv.e}. Finally, part \textbf{(f)} follows from Theorem
\ref{thm.sol.bininv.f} (since the maps $B$ and $W$ from Exercise
\ref{exe.bininv} \textbf{(f)} are restrictions of the maps $B$ and $W$ from
Definition \ref{def.sol.bininv.BW}). Hence, Exercise \ref{exe.bininv} is solved.
\end{proof}

\subsection{Solution to Exercise \ref{exe.binom.Hn-altsum}}

Before we solve Exercise \ref{exe.binom.Hn-altsum}, we state a lemma:

\begin{lemma}
\label{lem.sol.binom.Hn-altsum.1}Let $n\in\mathbb{N}$. Then,%
\[
\sum_{k=1}^{n+1}\dfrac{\left(  -1\right)  ^{k-1}}{k}\dbinom{n+1}{k}=\sum
_{k=1}^{n}\dfrac{\left(  -1\right)  ^{k-1}}{k}\dbinom{n}{k}+\dfrac{1}{n+1}.
\]

\end{lemma}

\begin{proof}
[Proof of Lemma \ref{lem.sol.binom.Hn-altsum.1}.]Clearly, $n+1$ is a positive
integer (since $n\in\mathbb{N}$). In other words, $n+1\in\left\{
1,2,3,\ldots\right\}  $. Hence, $n+1\neq0$. Thus, the fraction $\dfrac{1}%
{n+1}$ is well-defined.

Also, for each $k\in\left\{  1,2,3,\ldots\right\}  $, the fraction
$\dfrac{\left(  -1\right)  ^{k-1}}{k}$ is well-defined (since $k\neq0$).

Now, let $k\in\left\{  1,2,3,\ldots\right\}  $. Then, Proposition
\ref{prop.binom.X-1} (applied to $n+1$ and $k$ instead of $m$ and $n$) yields
$\dbinom{n+1}{k}=\dfrac{n+1}{k}\dbinom{\left(  n+1\right)  -1}{k-1}%
=\dfrac{n+1}{k}\dbinom{n}{k-1}$ (since $\left(  n+1\right)  -1=n$). Dividing
both sides of this equality by $n+1$, we find%
\begin{equation}
\dfrac{1}{n+1}\dbinom{n+1}{k}=\dfrac{1}{n+1}\cdot\dfrac{n+1}{k}\dbinom{n}%
{k-1}=\dfrac{1}{k}\dbinom{n}{k-1}. \label{sol.binom.Hn-altsum.1}%
\end{equation}
Now,
\begin{align}
\dfrac{\left(  -1\right)  ^{k-1}}{k}\dbinom{n}{k-1}  &  =\underbrace{\left(
-1\right)  ^{k-1}}_{=-\left(  -1\right)  ^{k}}\cdot\underbrace{\dfrac{1}%
{k}\dbinom{n}{k-1}}_{\substack{=\dfrac{1}{n+1}\dbinom{n+1}{k}\\\text{(by
(\ref{sol.binom.Hn-altsum.1}))}}}=-\left(  -1\right)  ^{k}\cdot\dfrac{1}%
{n+1}\dbinom{n+1}{k}\nonumber\\
&  =\dfrac{-1}{n+1}\left(  -1\right)  ^{k}\dbinom{n+1}{k}.
\label{sol.binom.Hn-altsum.2}%
\end{align}

Now, forget that we fixed $k$. We thus have proven the equality
(\ref{sol.binom.Hn-altsum.2}) for each $k\in\left\{  1,2,3,\ldots\right\}  $.

On the other hand, Proposition \ref{prop.binom.bin-id} \textbf{(c)} (applied
to $n+1$ instead of $n$) yields%
\[
\sum_{k=0}^{n+1}\left(  -1\right)  ^{k}\dbinom{n+1}{k}=%
\begin{cases}
1, & \text{if }n+1=0;\\
0, & \text{if }n+1\neq0
\end{cases}
=0
\]
(since $n+1\neq0$). Hence,%
\begin{align*}
0  &  =\sum_{k=0}^{n+1}\left(  -1\right)  ^{k}\dbinom{n+1}{k}%
=\underbrace{\left(  -1\right)  ^{0}}_{=1}\underbrace{\dbinom{n+1}{0}%
}_{\substack{=1\\\text{(by (\ref{eq.binom.00}) (applied to }m=n+1\text{))}%
}}+\sum_{k=1}^{n+1}\left(  -1\right)  ^{k}\dbinom{n+1}{k}\\
&  \ \ \ \ \ \ \ \ \ \ \left(  \text{here, we have split off the addend for
}k=0\text{ from the sum}\right) \\
&  =1+\sum_{k=1}^{n+1}\left(  -1\right)  ^{k}\dbinom{n+1}{k}.
\end{align*}
Subtracting $1$ from this equality, we obtain%
\begin{equation}
-1=\sum_{k=1}^{n+1}\left(  -1\right)  ^{k}\dbinom{n+1}{k}.
\label{sol.binom.Hn-altsum.4}%
\end{equation}

Also, $n<n+1$. Hence, Proposition \ref{prop.binom.0} (applied to $n$ and $n+1$
instead of $m$ and $n$) yields%
\[
\dbinom{n}{n+1}=0.
\]
Now,%
\begin{align}
\sum_{k=1}^{n+1}\dfrac{\left(  -1\right)  ^{k-1}}{k}\dbinom{n}{k}  &
=\sum_{k=1}^{n}\dfrac{\left(  -1\right)  ^{k-1}}{k}\dbinom{n}{k}%
+\dfrac{\left(  -1\right)  ^{\left(  n+1\right)  -1}}{n+1}\underbrace{\dbinom
{n}{n+1}}_{=0}\nonumber\\
&  \ \ \ \ \ \ \ \ \ \ \left(  \text{here, we have split off the addend for
}k=n+1\text{ from the sum}\right) \nonumber\\
&  =\sum_{k=1}^{n}\dfrac{\left(  -1\right)  ^{k-1}}{k}\dbinom{n}%
{k}+\underbrace{\dfrac{\left(  -1\right)  ^{\left(  n+1\right)  -1}}{n+1}%
0}_{=0}\nonumber\\
&  =\sum_{k=1}^{n}\dfrac{\left(  -1\right)  ^{k-1}}{k}\dbinom{n}{k}.
\label{sol.binom.Hn-altsum.6}%
\end{align}

Finally, each $k\in\left\{  1,2,3,\ldots\right\}  $ satisfies%
\begin{align}
\dbinom{n+1}{k}  &  =\dbinom{\left(  n+1\right)  -1}{k-1}+\dbinom{\left(
n+1\right)  -1}{k}\nonumber\\
&  \ \ \ \ \ \ \ \ \ \ \left(
\begin{array}
[c]{c}%
\text{by Proposition \ref{prop.binom.rec}}\\
\text{(applied to }n+1\text{ and }k\text{ instead of }m\text{ and }n\text{)}%
\end{array}
\right) \nonumber\\
&  =\dbinom{n}{k-1}+\dbinom{n}{k}\ \ \ \ \ \ \ \ \ \ \left(  \text{since
}\left(  n+1\right)  -1=n\right) \nonumber\\
&  =\dbinom{n}{k}+\dbinom{n}{k-1}. \label{sol.binom.Hn-altsum.8}%
\end{align}

Now,%
\begin{align*}
&  \sum_{k=1}^{n+1}\dfrac{\left(  -1\right)  ^{k-1}}{k}\underbrace{\dbinom
{n+1}{k}}_{\substack{=\dbinom{n}{k}+\dbinom{n}{k-1}\\\text{(by
(\ref{sol.binom.Hn-altsum.8}))}}}\\
&  =\sum_{k=1}^{n+1}\dfrac{\left(  -1\right)  ^{k-1}}{k}\left(  \dbinom{n}%
{k}+\dbinom{n}{k-1}\right)  =\underbrace{\sum_{k=1}^{n+1}\dfrac{\left(
-1\right)  ^{k-1}}{k}\dbinom{n}{k}}_{\substack{=\sum_{k=1}^{n}\dfrac{\left(
-1\right)  ^{k-1}}{k}\dbinom{n}{k}\\\text{(by (\ref{sol.binom.Hn-altsum.6}))}%
}}+\sum_{k=1}^{n+1}\underbrace{\dfrac{\left(  -1\right)  ^{k-1}}{k}\dbinom
{n}{k-1}}_{\substack{=\dfrac{-1}{n+1}\left(  -1\right)  ^{k}\dbinom{n+1}%
{k}\\\text{(by (\ref{sol.binom.Hn-altsum.2}))}}}\\
&  =\sum_{k=1}^{n}\dfrac{\left(  -1\right)  ^{k-1}}{k}\dbinom{n}%
{k}+\underbrace{\sum_{k=1}^{n+1}\dfrac{-1}{n+1}\left(  -1\right)  ^{k}%
\dbinom{n+1}{k}}_{=\dfrac{-1}{n+1}\sum_{k=1}^{n+1}\left(  -1\right)
^{k}\dbinom{n+1}{k}}\\
&  =\sum_{k=1}^{n}\dfrac{\left(  -1\right)  ^{k-1}}{k}\dbinom{n}{k}+\dfrac
{-1}{n+1}\underbrace{\sum_{k=1}^{n+1}\left(  -1\right)  ^{k}\dbinom{n+1}{k}%
}_{\substack{=-1\\\text{(by (\ref{sol.binom.Hn-altsum.4}))}}}\\
&  =\sum_{k=1}^{n}\dfrac{\left(  -1\right)  ^{k-1}}{k}\dbinom{n}{k}+\dfrac
{-1}{n+1}\left(  -1\right)  =\sum_{k=1}^{n}\dfrac{\left(  -1\right)  ^{k-1}%
}{k}\dbinom{n}{k}+\dfrac{1}{n+1}.
\end{align*}
This proves Lemma \ref{lem.sol.binom.Hn-altsum.1}.
\end{proof}

\begin{proof}
[Solution to Exercise \ref{exe.binom.Hn-altsum}.]We shall solve Exercise
\ref{exe.binom.Hn-altsum} by induction on $n$:

\textit{Induction base:} Comparing $\sum_{k=1}^{0}\dfrac{\left(  -1\right)
^{k-1}}{k}\dbinom{0}{k}=\left(  \text{empty sum}\right)  =0$ with $\dfrac
{1}{1}+\dfrac{1}{2}+\cdots+\dfrac{1}{0}=\left(  \text{empty sum}\right)  =0$,
we obtain $\sum_{k=1}^{0}\dfrac{\left(  -1\right)  ^{k-1}}{k}\dbinom{0}%
{k}=\dfrac{1}{1}+\dfrac{1}{2}+\cdots+\dfrac{1}{0}$. In other words, Exercise
\ref{exe.binom.Hn-altsum} holds for $n=0$. This completes the induction base.

\textit{Induction step:} Let $m\in\mathbb{N}$. Assume that Exercise
\ref{exe.binom.Hn-altsum} holds for $n=m$. We must prove that Exercise
\ref{exe.binom.Hn-altsum} holds for $n=m+1$.

We have assumed that Exercise \ref{exe.binom.Hn-altsum} holds for $n=m$. In
other words, we have%
\[
\sum_{k=1}^{m}\dfrac{\left(  -1\right)  ^{k-1}}{k}\dbinom{m}{k}=\dfrac{1}%
{1}+\dfrac{1}{2}+\cdots+\dfrac{1}{m}.
\]

Now, Lemma \ref{lem.sol.binom.Hn-altsum.1} (applied to $n=m$) yields%
\begin{align*}
\sum_{k=1}^{m+1}\dfrac{\left(  -1\right)  ^{k-1}}{k}\dbinom{m+1}{k}  &
=\underbrace{\sum_{k=1}^{m}\dfrac{\left(  -1\right)  ^{k-1}}{k}\dbinom{m}{k}%
}_{=\dfrac{1}{1}+\dfrac{1}{2}+\cdots+\dfrac{1}{m}}+\dfrac{1}{m+1}\\
&  =\left(  \dfrac{1}{1}+\dfrac{1}{2}+\cdots+\dfrac{1}{m}\right)  +\dfrac
{1}{m+1}=\dfrac{1}{1}+\dfrac{1}{2}+\cdots+\dfrac{1}{m+1}.
\end{align*}
In other words, Exercise \ref{exe.binom.Hn-altsum} holds for $n=m+1$. This
completes the induction step. Thus, Exercise \ref{exe.binom.Hn-altsum} is solved.
\end{proof}

\subsection{Solution to Exercise \ref{exe.binid.kurlis}}

We begin with a few lemmas before we come to the solution of Exercise
\ref{exe.binid.kurlis}. The first lemma is essentially trivial:

\begin{lemma}
\label{lem.sol.binid.kurlis.nz}Let $n\in\mathbb{N}$ and $k\in\left\{
0,1,\ldots,n\right\}  $. Then, $\dbinom{n}{k}$ is a positive integer.
\end{lemma}

\begin{proof}
[Proof of Lemma \ref{lem.sol.binid.kurlis.nz}.]From $n\in\mathbb{N}%
\subseteq\mathbb{Z}$ and $k\in\left\{  0,1,\ldots,n\right\}  \subseteq
\mathbb{N}$, we conclude that $\dbinom{n}{k}\in\mathbb{Z}$ (by Proposition
\ref{prop.binom.int} (applied to $n$ and $k$ instead of $m$ and $n$)). In
other words, $\dbinom{n}{k}$ is an integer.

From $k\in\left\{  0,1,\ldots,n\right\}  $, we obtain $k\leq n$, so that
$n\geq k$. Also, $k\in\left\{  0,1,\ldots,n\right\}  \subseteq\mathbb{N}$.
Hence, Proposition \ref{prop.binom.formula} (applied to $n$ and $k$ instead of
$m$ and $n$) yields $\dbinom{n}{k}=\dfrac{n!}{k!\left(  n-k\right)  !}$. But
$\dfrac{n!}{k!\left(  n-k\right)  !}$ is a positive rational number (since
$n!$, $k!$ and $\left(  n-k\right)  !$ are positive integers). In other words,
$\dbinom{n}{k}$ is a positive rational number (since $\dbinom{n}{k}=\dfrac
{n!}{k!\left(  n-k\right)  !}$). Hence, $\dbinom{n}{k}$ is a positive integer
(since $\dbinom{n}{k}$ is an integer). This proves Lemma
\ref{lem.sol.binid.kurlis.nz}.
\end{proof}

\begin{lemma}
\label{lem.sol.binid.kurlis.abs0}Let $q\in\mathbb{Q}$ and $m\in\mathbb{N}$ be
such that $q\neq m$. Then,%
\[
\dbinom{q}{m}=\dfrac{q}{q-m}\dbinom{q-1}{m}.
\]

\end{lemma}

\begin{proof}
[Proof of Lemma \ref{lem.sol.binid.kurlis.abs0}.]We have $q-m\neq0$ (since
$q\neq m$). Hence, the fraction $\dfrac{q}{q-m}$ is well-defined.

The equality (\ref{eq.binom.mn}) (applied to $q$ and $m$ instead of $m$ and
$n$) yields%
\begin{equation}
\dbinom{q}{m}=\dfrac{q\left(  q-1\right)  \cdots\left(  q-m+1\right)  }{m!}.
\label{pf.lem.sol.binid.kurlis.abs0.L}%
\end{equation}

The equality (\ref{eq.binom.mn}) (applied to $q-1$ and $m$ instead of $m$ and
$n$) yields%
\begin{align*}
\dbinom{q-1}{m}  &  =\dfrac{\left(  q-1\right)  \left(  \left(  q-1\right)
-1\right)  \cdots\left(  \left(  q-1\right)  -m+1\right)  }{m!}\\
&  =\dfrac{1}{m!}\cdot\underbrace{\left(  \left(  q-1\right)  \left(  \left(
q-1\right)  -1\right)  \cdots\left(  \left(  q-1\right)  -m+1\right)  \right)
}_{\substack{=\left(  q-1\right)  \left(  q-2\right)  \cdots\left(
q-m\right)  \\\text{(since }\left(  q-1\right)  -1=q-2\text{ and }\left(
q-1\right)  -m+1=q-m\text{)}}}\\
&  =\dfrac{1}{m!}\cdot\left(  \left(  q-1\right)  \left(  q-2\right)
\cdots\left(  q-m\right)  \right)  .
\end{align*}
Multiplying both sides of this equality by $\dfrac{q}{q-m}$, we find%
\begin{align*}
\dfrac{q}{q-m}\dbinom{q-1}{m}  &  =\dfrac{q}{q-m}\cdot\dfrac{1}{m!}%
\cdot\left(  \left(  q-1\right)  \left(  q-2\right)  \cdots\left(  q-m\right)
\right) \\
&  =\dfrac{1}{q-m}\cdot\dfrac{1}{m!}\cdot\underbrace{q\cdot\left(  \left(
q-1\right)  \left(  q-2\right)  \cdots\left(  q-m\right)  \right)
}_{\substack{=q\left(  q-1\right)  \cdots\left(  q-m\right)  \\=\left(
q\left(  q-1\right)  \cdots\left(  q-m+1\right)  \right)  \cdot\left(
q-m\right)  }}\\
&  =\dfrac{1}{q-m}\cdot\dfrac{1}{m!}\cdot\left(  q\left(  q-1\right)
\cdots\left(  q-m+1\right)  \right)  \cdot\left(  q-m\right) \\
&  =\dfrac{1}{m!}\cdot\left(  q\left(  q-1\right)  \cdots\left(  q-m+1\right)
\right)  =\dfrac{q\left(  q-1\right)  \cdots\left(  q-m+1\right)  }{m!}.
\end{align*}
Comparing this with (\ref{pf.lem.sol.binid.kurlis.abs0.L}), we obtain
$\dbinom{q}{m}=\dfrac{q}{q-m}\dbinom{q-1}{m}$. This proves Lemma
\ref{lem.sol.binid.kurlis.abs0}.
\end{proof}

\begin{lemma}
\label{lem.sol.binid.kurlis.abs}Let $n\in\mathbb{N}$ and $k\in\left\{
0,1,\ldots,n\right\}  $. Then:

\textbf{(a)} We have%
\[
\dbinom{n+1}{k+1}=\dfrac{n+1}{k+1}\dbinom{n}{k}.
\]

\textbf{(b)} We have
\[
\dbinom{n+1}{k}=\dfrac{n+1}{n+1-k}\dbinom{n}{k}.
\]

\end{lemma}

\begin{proof}
[Proof of Lemma \ref{lem.sol.binid.kurlis.abs}.]We have $k\in\left\{
0,1,\ldots,n\right\}  $, thus $k\leq n<n+1$. Hence, $n+1-k>0$. Thus, the
fraction $\dfrac{n+1}{n+1-k}$ is well-defined.

Also, from $k\in\left\{  0,1,\ldots,n\right\}  $, we obtain $k+1\in\left\{
1,2,\ldots,n+1\right\}  $; thus, $k+1>0$. Hence, the fraction $\dfrac
{n+1}{k+1}$ is well-defined.

We have $k\in\left\{  0,1,\ldots,n\right\}  $, so that $k+1\in\left\{
1,2,\ldots,n+1\right\}  \subseteq\left\{  1,2,3,\ldots\right\}  $. Thus,
Proposition \ref{prop.binom.X-1} (applied to $n+1$ and $k+1$ instead of $m$
and $n$) yields
\begin{equation}
\dbinom{n+1}{k+1}=\dfrac{n+1}{k+1}\dbinom{\left(  n+1\right)  -1}{\left(
k+1\right)  -1}=\dfrac{n+1}{k+1}\dbinom{n}{k}%
\end{equation}
(since $\left(  k+1\right)  -1=k$ and $\left(  n+1\right)  -1 = n$). This
proves Lemma \ref{lem.sol.binid.kurlis.abs} \textbf{(a)}.

\textbf{(b)} We have $k<n+1$ and thus $k\neq n+1$, so that $n+1\neq k$. Also,
$n+1\in\mathbb{N}\subseteq\mathbb{Q}$. Hence, Lemma
\ref{lem.sol.binid.kurlis.abs0} (applied to $q=n+1$ and $m=k$) yields%
\[
\dbinom{n+1}{k}=\dfrac{n+1}{n+1-k}\dbinom{n+1-1}{k}=\dfrac{n+1}{n+1-k}%
\dbinom{n}{k}%
\]
(since $n+1-1=n$). This proves Lemma \ref{lem.sol.binid.kurlis.abs}
\textbf{(b)}.
\end{proof}

Our next lemma is a simple consequence of the recurrence of the binomial coefficients:

\begin{lemma}
\label{lem.sol.binid.kurlis.rec}Let $n\in\mathbb{N}$. Let $k\in\left\{
0,1,\ldots,n\right\}  $. Then,%
\begin{equation}
\dfrac{1}{\dbinom{n}{k}}=\left(  \dfrac{1}{\dbinom{n+1}{k}}+\dfrac{1}%
{\dbinom{n+1}{k+1}}\right)  \dfrac{n+1}{n+2}.
\label{eq.lem.sol.binid.kurlis.rec.eq}%
\end{equation}
(In particular, all the fractions $\dfrac{1}{\dbinom{n}{k}}$, $\dfrac
{1}{\dbinom{n+1}{k}}$ and $\dfrac{1}{\dbinom{n+1}{k+1}}$ in this equality are well-defined.)
\end{lemma}

\begin{proof}
[Proof of Lemma \ref{lem.sol.binid.kurlis.rec}.]Lemma
\ref{lem.sol.binid.kurlis.nz} shows that $\dbinom{n}{k}$ is a positive
integer. Thus, $\dbinom{n}{k}\neq0$. Hence, the fraction $\dfrac{1}{\dbinom
{n}{k}}$ is well-defined.

We have $k\in\left\{  0,1,\ldots,n\right\}  \subseteq\left\{  0,1,\ldots
,n+1\right\}  $. Thus, Lemma \ref{lem.sol.binid.kurlis.nz} (applied to $n+1$
instead of $n$) shows that $\dbinom{n+1}{k}$ is a positive integer. Thus,
$\dbinom{n+1}{k}\neq0$. Hence, the fraction $\dfrac{1}{\dbinom{n+1}{k}}$ is well-defined.

We have $k\in\left\{  0,1,\ldots,n\right\}  $, thus $k+1\in\left\{
1,2,\ldots,n+1\right\}  \subseteq\left\{  0,1,\ldots,n+1\right\}  $. Thus,
Lemma \ref{lem.sol.binid.kurlis.nz} (applied to $n+1$ and $k+1$ instead of $n$
and $k$) shows that $\dbinom{n+1}{k+1}$ is a positive integer. Thus,
$\dbinom{n+1}{k+1}\neq0$. Hence, the fraction $\dfrac{1}{\dbinom{n+1}{k+1}}$
is well-defined.

It remains to prove the equality (\ref{eq.lem.sol.binid.kurlis.rec.eq}). We
have%
\begin{align*}
&  \dfrac{1}{\dbinom{n+1}{k}}+\dfrac{1}{\dbinom{n+1}{k+1}}\\
&  =1/\underbrace{\dbinom{n+1}{k}}_{\substack{=\dfrac{n+1}{n+1-k}\dbinom{n}%
{k}\\\text{(by Lemma \ref{lem.sol.binid.kurlis.abs} \textbf{(b)})}%
}}+1/\underbrace{\dbinom{n+1}{k+1}}_{\substack{=\dfrac{n+1}{k+1}\dbinom{n}%
{k}\\\text{(by Lemma \ref{lem.sol.binid.kurlis.abs} \textbf{(a)})}}}\\
&  =\underbrace{1/\left(  \dfrac{n+1}{n+1-k}\dbinom{n}{k}\right)  }%
_{=\dfrac{n+1-k}{n+1}/\dbinom{n}{k}}+\underbrace{1/\left(  \dfrac{n+1}%
{k+1}\dbinom{n}{k}\right)  }_{=\dfrac{k+1}{n+1}/\dbinom{n}{k}}\\
&  =\dfrac{n+1-k}{n+1}/\dbinom{n}{k}+\dfrac{k+1}{n+1}/\dbinom{n}%
{k}=\underbrace{\left(  \dfrac{n+1-k}{n+1}+\dfrac{k+1}{n+1}\right)  }%
_{=\dfrac{n+2}{n+1}}/\dbinom{n}{k}\\
&  =\dfrac{n+2}{n+1}/\dbinom{n}{k}.
\end{align*}
Multiplying both sides of this equality by $\dfrac{n+1}{n+2}$, we obtain%
\[
\left(  \dfrac{1}{\dbinom{n+1}{k}}+\dfrac{1}{\dbinom{n+1}{k+1}}\right)
\dfrac{n+1}{n+2}=\left(  \dfrac{n+2}{n+1}/\dbinom{n}{k}\right)  \dfrac
{n+1}{n+2}=\dfrac{1}{\dbinom{n}{k}}.
\]
This proves Lemma \ref{lem.sol.binid.kurlis.rec}.
\end{proof}

We can now prove part \textbf{(a)} of Exercise \ref{exe.binid.kurlis}:

\begin{proposition}
\label{prop.sol.binid.kurlis.a}Let $n\in\mathbb{N}$.

\textbf{(a)} We have%
\[
\sum_{k=0}^{n}\dfrac{\left(  -1\right)  ^{k}}{\dbinom{n}{k}}=\dfrac{n+1}%
{n+2}\left(  1+\left(  -1\right)  ^{n}\right)  .
\]

\textbf{(b)} We have
\[
\sum_{k=0}^{n}\dfrac{\left(  -1\right)  ^{k}}{\dbinom{n}{k}}=2\cdot\dfrac
{n+1}{n+2}\left[  n\text{ is even}\right]  .
\]
(Here, we are using the Iverson bracket notation, as in Definition
\ref{def.iverson}; thus, $\left[  n\text{ is even}\right]  $ is $1$ if $n$ is
even and $0$ otherwise.)
\end{proposition}

\begin{proof}
[Proof of Proposition \ref{prop.sol.binid.kurlis.a}.]\textbf{(a)} We have
$\dbinom{n+1}{0}=1$ (by Proposition \ref{prop.binom.00} \textbf{(a)} (applied
to $m=n+1$)), thus $\dfrac{1}{\dbinom{n+1}{0}}=\dfrac{1}{1}=1$.

Also, $\dbinom{n+1}{n+1}=1$ (by Proposition \ref{prop.binom.mm} (applied to
$m=n+1$)), thus $\dfrac{1}{\dbinom{n+1}{n+1}}=\dfrac{1}{1}=1$.

We have%
\begin{align}
&  \sum_{k=0}^{n}\underbrace{\dfrac{\left(  -1\right)  ^{k}}{\dbinom{n}{k}}%
}_{=\left(  -1\right)  ^{k}\cdot\dfrac{1}{\dbinom{n}{k}}}\nonumber\\
&  =\sum_{k=0}^{n}\left(  -1\right)  ^{k}\cdot\underbrace{\dfrac{1}{\dbinom
{n}{k}}}_{\substack{=\left(  \dfrac{1}{\dbinom{n+1}{k}}+\dfrac{1}{\dbinom
{n+1}{k+1}}\right)  \dfrac{n+1}{n+2}\\\text{(by
(\ref{eq.lem.sol.binid.kurlis.rec.eq}))}}}\nonumber\\
&  =\sum_{k=0}^{n}\left(  -1\right)  ^{k}\cdot\left(  \dfrac{1}{\dbinom
{n+1}{k}}+\dfrac{1}{\dbinom{n+1}{k+1}}\right)  \dfrac{n+1}{n+2}\nonumber\\
&  =\dfrac{n+1}{n+2}\sum_{k=0}^{n}\left(  -1\right)  ^{k}\cdot\left(
\dfrac{1}{\dbinom{n+1}{k}}+\dfrac{1}{\dbinom{n+1}{k+1}}\right)  .
\label{pf.prop.sol.binid.kurlis.a.1}%
\end{align}

But%
\begin{align*}
&  \sum_{k=0}^{n}\underbrace{\left(  -1\right)  ^{k}\cdot\left(  \dfrac
{1}{\dbinom{n+1}{k}}+\dfrac{1}{\dbinom{n+1}{k+1}}\right)  }_{=\left(
-1\right)  ^{k}\cdot\dfrac{1}{\dbinom{n+1}{k}}+\left(  -1\right)  ^{k}%
\cdot\dfrac{1}{\dbinom{n+1}{k+1}}}\\
&  =\sum_{k=0}^{n}\left(  \left(  -1\right)  ^{k}\cdot\dfrac{1}{\dbinom
{n+1}{k}}+\left(  -1\right)  ^{k}\cdot\dfrac{1}{\dbinom{n+1}{k+1}}\right) \\
&  =\underbrace{\sum_{k=0}^{n}\left(  -1\right)  ^{k}\cdot\dfrac{1}%
{\dbinom{n+1}{k}}}_{\substack{=\left(  -1\right)  ^{0}\cdot\dfrac{1}%
{\dbinom{n+1}{0}}+\sum_{k=1}^{n}\left(  -1\right)  ^{k}\cdot\dfrac{1}%
{\dbinom{n+1}{k}}\\\text{(here, we have split off the addend}\\\text{for
}k=0\text{ from the sum)}}}+\underbrace{\sum_{k=0}^{n}\left(  -1\right)
^{k}\cdot\dfrac{1}{\dbinom{n+1}{k+1}}}_{\substack{=\left(  -1\right)
^{n}\cdot\dfrac{1}{\dbinom{n+1}{n+1}}+\sum_{k=0}^{n-1}\left(  -1\right)
^{k}\cdot\dfrac{1}{\dbinom{n+1}{k+1}}\\\text{(here, we have split off the
addend}\\\text{for }k=n\text{ from the sum)}}}\\
&  =\underbrace{\left(  -1\right)  ^{0}}_{=1}\cdot\underbrace{\dfrac
{1}{\dbinom{n+1}{0}}}_{=1}+\underbrace{\sum_{k=1}^{n}\left(  -1\right)
^{k}\cdot\dfrac{1}{\dbinom{n+1}{k}}}_{\substack{=\sum_{k=0}^{n-1}\left(
-1\right)  ^{k+1}\cdot\dfrac{1}{\dbinom{n+1}{k+1}}\\\text{(here, we have
substituted }k+1\\\text{for }k\text{ in the sum)}}}+\left(  -1\right)
^{n}\cdot\underbrace{\dfrac{1}{\dbinom{n+1}{n+1}}}_{=1}+\sum_{k=0}%
^{n-1}\left(  -1\right)  ^{k}\cdot\dfrac{1}{\dbinom{n+1}{k+1}}\\
&  =1+\sum_{k=0}^{n-1}\underbrace{\left(  -1\right)  ^{k+1}}_{=-\left(
-1\right)  ^{k}}\cdot\dfrac{1}{\dbinom{n+1}{k+1}}+\left(  -1\right)  ^{n}%
+\sum_{k=0}^{n-1}\left(  -1\right)  ^{k}\cdot\dfrac{1}{\dbinom{n+1}{k+1}}\\
&  =1+\underbrace{\sum_{k=0}^{n-1}\left(  -\left(  -1\right)  ^{k}\right)
\cdot\dfrac{1}{\dbinom{n+1}{k+1}}}_{=-\sum_{k=0}^{n-1}\left(  -1\right)
^{k}\cdot\dfrac{1}{\dbinom{n+1}{k+1}}}+\left(  -1\right)  ^{n}+\sum
_{k=0}^{n-1}\left(  -1\right)  ^{k}\cdot\dfrac{1}{\dbinom{n+1}{k+1}}\\
&  =1+\left(  -\sum_{k=0}^{n-1}\left(  -1\right)  ^{k}\cdot\dfrac{1}%
{\dbinom{n+1}{k+1}}\right)  +\left(  -1\right)  ^{n}+\sum_{k=0}^{n-1}\left(
-1\right)  ^{k}\cdot\dfrac{1}{\dbinom{n+1}{k+1}}=1+\left(  -1\right)  ^{n}.
\end{align*}
Hence, (\ref{pf.prop.sol.binid.kurlis.a.1}) becomes%
\begin{align*}
\sum_{k=0}^{n}\dfrac{\left(  -1\right)  ^{k}}{\dbinom{n}{k}}  &  =\dfrac
{n+1}{n+2}\underbrace{\sum_{k=0}^{n}\left(  -1\right)  ^{k}\cdot\left(
\dfrac{1}{\dbinom{n+1}{k}}+\dfrac{1}{\dbinom{n+1}{k+1}}\right)  }_{=1+\left(
-1\right)  ^{n}}\\
&  =\dfrac{n+1}{n+2}\left(  1+\left(  -1\right)  ^{n}\right)  .
\end{align*}
This proves Proposition \ref{prop.sol.binid.kurlis.a} \textbf{(a)}.

\textbf{(b)} The definition of the truth value $\left[  n\text{ is
even}\right]  $ shows that%
\begin{equation}
\left[  n\text{ is even}\right]  =%
\begin{cases}
1, & \text{if }n\text{ is even;}\\
0, & \text{if }n\text{ is not even}%
\end{cases}
=%
\begin{cases}
1, & \text{if }n\text{ is even;}\\
0, & \text{if }n\text{ is odd}%
\end{cases}
\label{pf.prop.sol.binid.kurlis.a.iverson}%
\end{equation}
(because the condition \textquotedblleft$n$ is not even\textquotedblright\ is
equivalent to \textquotedblleft$n$ is odd\textquotedblright). But Proposition
\ref{prop.sol.binid.kurlis.a} \textbf{(a)} yields%
\begin{align*}
\sum_{k=0}^{n}\dfrac{\left(  -1\right)  ^{k}}{\dbinom{n}{k}}  &  =\dfrac
{n+1}{n+2}\left(  1+\left(  -1\right)  ^{n}\right)  =2\cdot\dfrac{n+1}%
{n+2}\cdot\underbrace{\dfrac{1}{2}\left(  1+\left(  -1\right)  ^{n}\right)
}_{\substack{=%
\begin{cases}
1, & \text{if }n\text{ is even;}\\
0, & \text{if }n\text{ is odd}%
\end{cases}
\\\text{(by Lemma \ref{lem.sol.bininv.pendulum} \textbf{(a)})}}}\\
&  =2\cdot\dfrac{n+1}{n+2}\cdot\underbrace{%
\begin{cases}
1, & \text{if }n\text{ is even;}\\
0, & \text{if }n\text{ is odd}%
\end{cases}
}_{\substack{=\left[  n\text{ is even}\right]  \\\text{(by
(\ref{pf.prop.sol.binid.kurlis.a.iverson}))}}}=2\cdot\dfrac{n+1}{n+2}\left[
n\text{ is even}\right]  .
\end{align*}
This proves Proposition \ref{prop.sol.binid.kurlis.a} \textbf{(b)}.
\end{proof}

We remark that Proposition \ref{prop.sol.binid.kurlis.a} \textbf{(a)} can also
be seen as a particular case of the identity%
\[
\sum_{k=0}^{n}\dfrac{\left(  -1\right)  ^{k}}{\dbinom{x}{k}}=\left(
1+\dfrac{\left(  -1\right)  ^{n}}{\dbinom{x+1}{n+1}}\right)  \dfrac{x+1}{x+2}%
\]
(for all $n\in\left\{  -1,0,1,\ldots\right\}  $ and $x\in\mathbb{Q}%
\setminus\left\{  -2,-1,0,\ldots,n-1\right\}  $), which was observed by user
\textquotedblleft user90369\textquotedblright\ on
\url{https://math.stackexchange.com/a/3251880/} and can be proven fairly
easily by induction on $n$.

\begin{proof}
[Solution to Exercise \ref{exe.binid.kurlis}.]\textbf{(a)} Proposition
\ref{prop.sol.binid.kurlis.a} \textbf{(b)} yields
\[
\sum_{k=0}^{n}\dfrac{\left(  -1\right)  ^{k}}{\dbinom{n}{k}}=2\cdot\dfrac
{n+1}{n+2}\left[  n\text{ is even}\right]  .
\]
This solves Exercise \ref{exe.binid.kurlis} \textbf{(a)}.

\textbf{(b)} Let us forget that we fixed $n$. We shall solve Exercise
\ref{exe.binid.kurlis} \textbf{(b)} by induction on $n$:

\textit{Induction base:} Proposition \ref{prop.binom.00} \textbf{(a)} (applied
to $m=0$) yields $\dbinom{0}{0}=1$. We have $0+1=1$ and thus%
\[
\dfrac{0+1}{2^{0+1}}\sum_{k=1}^{0+1}\dfrac{2^{k}}{k}=\dfrac{1}{2^{1}%
}\underbrace{\sum_{k=1}^{1}\dfrac{2^{k}}{k}}_{=\dfrac{2^{1}}{1}}=\dfrac
{1}{2^{1}}\cdot\dfrac{2^{1}}{1}=1.
\]
Comparing this with%
\[
\sum_{k=0}^{0}\dfrac{1}{\dbinom{0}{k}}=\dfrac{1}{\dbinom{0}{0}}%
=1/\underbrace{\dbinom{0}{0}}_{=1}=1/1=1,
\]
we obtain $\sum_{k=0}^{0}\dfrac{1}{\dbinom{0}{k}}=\dfrac{0+1}{2^{0+1}}%
\sum_{k=1}^{0+1}\dfrac{2^{k}}{k}$. In other words, Exercise
\ref{exe.binid.kurlis} \textbf{(b)} holds for $n=0$. This completes the
induction base.

\textit{Induction step:} Let $m\in\mathbb{N}$. Assume that Exercise
\ref{exe.binid.kurlis} \textbf{(b)} holds for $n=m$. We must prove that
Exercise \ref{exe.binid.kurlis} \textbf{(b)} holds for $n=m+1$.

We have assumed that Exercise \ref{exe.binid.kurlis} \textbf{(b)} holds for
$n=m$. In other words, we have%
\begin{equation}
\sum_{k=0}^{m}\dfrac{1}{\dbinom{m}{k}}=\dfrac{m+1}{2^{m+1}}\sum_{k=1}%
^{m+1}\dfrac{2^{k}}{k}. \label{sol.binid.kurlis.b.IH}%
\end{equation}

We have $\dbinom{m+1}{0}=1$ (by Proposition \ref{prop.binom.00} \textbf{(a)}
(applied to $m+1$ instead of $m$)), thus $\dfrac{1}{\dbinom{m+1}{0}}=\dfrac
{1}{1}=1$.

Also, $\dbinom{m+1}{m+1}=1$ (by Proposition \ref{prop.binom.mm} (applied to
$m+1$ instead of $m$)), thus $\dfrac{1}{\dbinom{m+1}{m+1}}=\dfrac{1}{1}=1$.

Notice that $m+1>m\geq0$ (since $m\in\mathbb{N}$), and thus $m+1\neq0$.

We have%
\begin{align*}
\sum_{k=0}^{m+1}\dfrac{1}{\dbinom{m+1}{k}}  &  =\sum_{k=0}^{m}\dfrac
{1}{\dbinom{m+1}{k}}+\underbrace{\dfrac{1}{\dbinom{m+1}{m+1}}}_{=1}\\
&  \ \ \ \ \ \ \ \ \ \ \left(
\begin{array}
[c]{c}%
\text{here, we have split off the addend for }k=m+1\\
\text{from the sum}%
\end{array}
\right) \\
&  =\sum_{k=0}^{m}\dfrac{1}{\dbinom{m+1}{k}}+1.
\end{align*}
Thus,%
\begin{equation}
\sum_{k=0}^{m}\dfrac{1}{\dbinom{m+1}{k}}=\sum_{k=0}^{m+1}\dfrac{1}%
{\dbinom{m+1}{k}}-1. \label{sol.binid.kurlis.b.sum1}%
\end{equation}

Also,%
\begin{align*}
\sum_{k=0}^{m+1}\dfrac{1}{\dbinom{m+1}{k}}  &  =\underbrace{\dfrac{1}%
{\dbinom{m+1}{0}}}_{=1}+\sum_{k=1}^{m+1}\dfrac{1}{\dbinom{m+1}{k}}\\
&  \ \ \ \ \ \ \ \ \ \ \left(
\begin{array}
[c]{c}%
\text{here, we have split off the addend for }k=0\\
\text{from the sum}%
\end{array}
\right) \\
&  =1+\sum_{k=1}^{m+1}\dfrac{1}{\dbinom{m+1}{k}}=1+\sum_{k=0}^{m}\dfrac
{1}{\dbinom{m+1}{k+1}}%
\end{align*}
(here, we have substituted $k+1$ for $k$ in the sum). Thus,%
\begin{equation}
\sum_{k=0}^{m}\dfrac{1}{\dbinom{m+1}{k+1}}=\sum_{k=0}^{m+1}\dfrac{1}%
{\dbinom{m+1}{k}}-1. \label{sol.binid.kurlis.b.sum2}%
\end{equation}

Every $k\in\left\{  0,1,\ldots,m\right\}  $ satisfies%
\begin{equation}
\dfrac{1}{\dbinom{m}{k}}=\left(  \dfrac{1}{\dbinom{m+1}{k}}+\dfrac{1}%
{\dbinom{m+1}{k+1}}\right)  \dfrac{m+1}{m+2} \label{sol.binid.kurlis.b.lem}%
\end{equation}
(by (\ref{eq.lem.sol.binid.kurlis.rec.eq}) (applied to $n=m$)).

Now,%
\begin{align*}
\sum_{k=0}^{m}\underbrace{\dfrac{1}{\dbinom{m}{k}}}_{=\left(  \dfrac
{1}{\dbinom{m+1}{k}}+\dfrac{1}{\dbinom{m+1}{k+1}}\right)  \dfrac{m+1}{m+2}}
&  =\sum_{k=0}^{m}\left(  \dfrac{1}{\dbinom{m+1}{k}}+\dfrac{1}{\dbinom
{m+1}{k+1}}\right)  \dfrac{m+1}{m+2}\\
&  =\dfrac{m+1}{m+2}\sum_{k=0}^{m}\left(  \dfrac{1}{\dbinom{m+1}{k}}+\dfrac
{1}{\dbinom{m+1}{k+1}}\right)  .
\end{align*}
In view of%
\begin{align*}
\sum_{k=0}^{m}\left(  \dfrac{1}{\dbinom{m+1}{k}}+\dfrac{1}{\dbinom{m+1}{k+1}%
}\right)   &  =\underbrace{\sum_{k=0}^{m}\dfrac{1}{\dbinom{m+1}{k}}%
}_{\substack{=\sum_{k=0}^{m+1}\dfrac{1}{\dbinom{m+1}{k}}-1\\\text{(by
(\ref{sol.binid.kurlis.b.sum1}))}}}+\underbrace{\sum_{k=0}^{m}\dfrac
{1}{\dbinom{m+1}{k+1}}}_{\substack{=\sum_{k=0}^{m+1}\dfrac{1}{\dbinom{m+1}{k}%
}-1\\\text{(by (\ref{sol.binid.kurlis.b.sum2}))}}}\\
&  =\left(  \sum_{k=0}^{m+1}\dfrac{1}{\dbinom{m+1}{k}}-1\right)  +\left(
\sum_{k=0}^{m+1}\dfrac{1}{\dbinom{m+1}{k}}-1\right) \\
&  =2\left(  \sum_{k=0}^{m+1}\dfrac{1}{\dbinom{m+1}{k}}-1\right)  ,
\end{align*}
this becomes%
\[
\sum_{k=0}^{m}\dfrac{1}{\dbinom{m}{k}}=\dfrac{m+1}{m+2}\underbrace{\sum
_{k=0}^{m}\left(  \dfrac{1}{\dbinom{m+1}{k}}+\dfrac{1}{\dbinom{m+1}{k+1}%
}\right)  }_{=2\left(  \sum_{k=0}^{m+1}\dfrac{1}{\dbinom{m+1}{k}}-1\right)
}=\dfrac{m+1}{m+2}\cdot2\left(  \sum_{k=0}^{m+1}\dfrac{1}{\dbinom{m+1}{k}%
}-1\right)  .
\]
Comparing this with (\ref{sol.binid.kurlis.b.IH}), we obtain%
\[
\dfrac{m+1}{m+2}\cdot2\left(  \sum_{k=0}^{m+1}\dfrac{1}{\dbinom{m+1}{k}%
}-1\right)  =\dfrac{m+1}{2^{m+1}}\sum_{k=1}^{m+1}\dfrac{2^{k}}{k}.
\]
Multiplying both sides of this equality by $\dfrac{m+2}{m+1}$ (this is
allowed, since $m+1\neq0$), we obtain%
\[
\dfrac{m+2}{m+1}\cdot\dfrac{m+1}{m+2}\cdot2\left(  \sum_{k=0}^{m+1}\dfrac
{1}{\dbinom{m+1}{k}}-1\right)  =\dfrac{m+2}{m+1}\cdot\dfrac{m+1}{2^{m+1}}%
\sum_{k=1}^{m+1}\dfrac{2^{k}}{k}=\dfrac{m+2}{2^{m+1}}\sum_{k=1}^{m+1}%
\dfrac{2^{k}}{k}.
\]
Hence,%
\[
\dfrac{m+2}{2^{m+1}}\sum_{k=1}^{m+1}\dfrac{2^{k}}{k}=\dfrac{m+2}{m+1}%
\cdot\dfrac{m+1}{m+2}\cdot2\left(  \sum_{k=0}^{m+1}\dfrac{1}{\dbinom{m+1}{k}%
}-1\right)  =2\left(  \sum_{k=0}^{m+1}\dfrac{1}{\dbinom{m+1}{k}}-1\right)  .
\]
Dividing both sides of this equality by $2$, we find%
\[
\dfrac{1}{2}\cdot\dfrac{m+2}{2^{m+1}}\sum_{k=1}^{m+1}\dfrac{2^{k}}{k}%
=\sum_{k=0}^{m+1}\dfrac{1}{\dbinom{m+1}{k}}-1.
\]
Solving this equality for $\sum_{k=0}^{m+1}\dfrac{1}{\dbinom{m+1}{k}}$, we
obtain%
\[
\sum_{k=0}^{m+1}\dfrac{1}{\dbinom{m+1}{k}}=\underbrace{\dfrac{1}{2}\cdot
\dfrac{m+2}{2^{m+1}}}_{\substack{=\dfrac{m+2}{2\cdot2^{m+1}}=\dfrac
{m+2}{2^{m+2}}\\\text{(since }2\cdot2^{m+1}=2^{\left(  m+1\right)  +1}%
=2^{m+2}\\\text{(because }\left(  m+1\right)  +1=m+2\text{))}}}\sum
_{k=1}^{m+1}\dfrac{2^{k}}{k}+1=\dfrac{m+2}{2^{m+2}}\sum_{k=1}^{m+1}%
\dfrac{2^{k}}{k}+1.
\]
Comparing this with%
\begin{align*}
&  \underbrace{\dfrac{\left(  m+1\right)  +1}{2^{\left(  m+1\right)  +1}}%
}_{\substack{=\dfrac{m+2}{2^{m+2}}\\\text{(since }\left(  m+1\right)
+1=m+2\text{)}}}\ \ \ \underbrace{\sum_{k=1}^{\left(  m+1\right)  +1}%
\dfrac{2^{k}}{k}}_{\substack{=\sum_{k=1}^{m+2}\dfrac{2^{k}}{k}\\\text{(since
}\left(  m+1\right)  +1=m+2\text{)}}}\\
&  =\dfrac{m+2}{2^{m+2}}\underbrace{\sum_{k=1}^{m+2}\dfrac{2^{k}}{k}%
}_{\substack{=\sum_{k=1}^{m+1}\dfrac{2^{k}}{k}+\dfrac{2^{m+2}}{m+2}%
\\\text{(here, we have split off the}\\\text{addend for }k=m+2\text{ from the
sum)}}}=\dfrac{m+2}{2^{m+2}}\left(  \sum_{k=1}^{m+1}\dfrac{2^{k}}{k}%
+\dfrac{2^{m+2}}{m+2}\right) \\
&  =\dfrac{m+2}{2^{m+2}}\sum_{k=1}^{m+1}\dfrac{2^{k}}{k}+\underbrace{\dfrac
{m+2}{2^{m+2}}\cdot\dfrac{2^{m+2}}{m+2}}_{=1}=\dfrac{m+2}{2^{m+2}}\sum
_{k=1}^{m+1}\dfrac{2^{k}}{k}+1,
\end{align*}
we obtain%
\[
\sum_{k=0}^{m+1}\dfrac{1}{\dbinom{m+1}{k}}=\dfrac{\left(  m+1\right)
+1}{2^{\left(  m+1\right)  +1}}\sum_{k=1}^{\left(  m+1\right)  +1}\dfrac
{2^{k}}{k}.
\]
In other words, Exercise \ref{exe.binid.kurlis} \textbf{(b)} holds for
$n=m+1$. This completes the induction step. Thus, Exercise
\ref{exe.binid.kurlis} \textbf{(b)} is solved by induction.
\end{proof}

\subsection{Solution to Exercise \ref{exe.AoPS333199}}

We shall prove the following generalization of Exercise \ref{exe.AoPS333199}:

\begin{proposition}
\label{prop.AoPS333199.gen}Let $\mathbb{K}$ be a commutative ring. (See
Definition \ref{def.commring} for the definition of a \textquotedblleft
commutative ring\textquotedblright. For example, we can set $\mathbb{K}%
=\mathbb{Z}$ or $\mathbb{K}=\mathbb{R}$ or $\mathbb{K}=\mathbb{Q}\left[
X\right]  $.) Let $x$ and $y$ be two elements of $\mathbb{K}$. For any
$n\in\mathbb{N}$ and $m\in\mathbb{N}$, define $Y_{m,n}\in\mathbb{K}$ by%
\[
Y_{m,n}=\sum_{k=0}^{n}y^{k}\dbinom{n}{k}\left(  x^{n-k}+y\right)  ^{m}.
\]
Then, $Y_{m,n}=Y_{n,m}$ for any $n\in\mathbb{N}$ and $m\in\mathbb{N}$.
\end{proposition}

\begin{proof}
[Proof of Proposition \ref{prop.AoPS333199.gen}.]For any $n\in\mathbb{N}$ and
$m\in\mathbb{N}$, we have%
\begin{align}
Y_{m,n}  &  =\sum_{k=0}^{n}y^{k}\dbinom{n}{k}\left(  x^{n-k}+y\right)
^{m}=\sum_{\ell=0}^{n}y^{\ell}\dbinom{n}{\ell}\left(  \underbrace{x^{n-\ell
}+y}_{=y+x^{n-\ell}}\right)  ^{m}\nonumber\\
&  \ \ \ \ \ \ \ \ \ \ \left(  \text{here, we have renamed the summation index
}k\text{ as }\ell\right) \nonumber\\
&  =\sum_{\ell=0}^{n}y^{\ell}\dbinom{n}{\ell}\underbrace{\left(  y+x^{n-\ell
}\right)  ^{m}}_{\substack{=\sum_{k=0}^{m}\dbinom{m}{k}y^{k}\left(  x^{n-\ell
}\right)  ^{m-k}\\\text{(by the binomial formula)}}}=\sum_{\ell=0}^{n}y^{\ell
}\dbinom{n}{\ell}\left(  \sum_{k=0}^{m}\dbinom{m}{k}y^{k}\left(  x^{n-\ell
}\right)  ^{m-k}\right) \nonumber\\
&  =\sum_{\ell=0}^{n}\sum_{k=0}^{m}\underbrace{y^{\ell}\dbinom{n}{\ell}%
\dbinom{m}{k}}_{=\dbinom{n}{\ell}\dbinom{m}{k}y^{\ell}}y^{k}\left(  x^{n-\ell
}\right)  ^{m-k}=\sum_{\ell=0}^{n}\sum_{k=0}^{m}\dbinom{n}{\ell}\dbinom{m}%
{k}\underbrace{y^{\ell}y^{k}}_{=y^{\ell+k}}\underbrace{\left(  x^{n-\ell
}\right)  ^{m-k}}_{=x^{\left(  n-\ell\right)  \left(  m-k\right)  }%
}\nonumber\\
&  =\sum_{\ell=0}^{n}\sum_{k=0}^{m}\dbinom{n}{\ell}\dbinom{m}{k}y^{\ell
+k}x^{\left(  n-\ell\right)  \left(  m-k\right)  }.
\label{pf.prop.AoPS333199.gen.1}%
\end{align}
Now, let $n\in\mathbb{N}$ and $m\in\mathbb{N}$. Then,
(\ref{pf.prop.AoPS333199.gen.1}) (applied to $m$ and $n$ instead of $n$ and
$m$) shows that%
\begin{align*}
Y_{n,m}  &  =\underbrace{\sum_{\ell=0}^{m}\sum_{k=0}^{n}}_{=\sum_{k=0}^{n}%
\sum_{\ell=0}^{m}}\ \ \underbrace{\dbinom{m}{\ell}\dbinom{n}{k}}_{=\dbinom
{n}{k}\dbinom{m}{\ell}}\ \ \underbrace{y^{\ell+k}}_{=y^{k+\ell}}%
\ \ \underbrace{x^{\left(  m-\ell\right)  \left(  n-k\right)  }}_{=x^{\left(
n-k\right)  \left(  m-\ell\right)  }}\\
&  =\sum_{k=0}^{n}\sum_{\ell=0}^{m}\dbinom{n}{k}\dbinom{m}{\ell}y^{k+\ell
}x^{\left(  n-k\right)  \left(  m-\ell\right)  }=\sum_{k=0}^{n}\sum_{g=0}%
^{m}\dbinom{n}{k}\dbinom{m}{g}y^{k+g}x^{\left(  n-k\right)  \left(
m-g\right)  }\\
&  \ \ \ \ \ \ \ \ \ \ \left(  \text{here, we have renamed the summation index
}\ell\text{ as }g\right) \\
&  =\sum_{\ell=0}^{n}\sum_{g=0}^{m}\dbinom{n}{\ell}\dbinom{m}{g}y^{\ell
+g}x^{\left(  n-\ell\right)  \left(  m-g\right)  }\\
&  \ \ \ \ \ \ \ \ \ \ \left(  \text{here, we have renamed the summation index
}k\text{ as }\ell\right) \\
&  =\sum_{\ell=0}^{n}\sum_{k=0}^{m}\dbinom{n}{\ell}\dbinom{m}{k}y^{\ell
+k}x^{\left(  n-\ell\right)  \left(  m-k\right)  }\\
&  \ \ \ \ \ \ \ \ \ \ \left(  \text{here, we have renamed the summation index
}g\text{ as }k\right) \\
&  =Y_{m,n}\ \ \ \ \ \ \ \ \ \ \left(  \text{by
(\ref{pf.prop.AoPS333199.gen.1})}\right)  .
\end{align*}
This proves Proposition \ref{prop.AoPS333199.gen}.
\end{proof}

\begin{proof}
[Solution to Exercise \ref{exe.AoPS333199}.]Set $\mathbb{K}=\mathbb{Z}\left[
X\right]  $, and define two elements $x$ and $y$ of $\mathbb{K}$ by $x=X$ and
$y=-1$. For any $n\in\mathbb{N}$ and $m\in\mathbb{N}$, define $Y_{m,n}%
\in\mathbb{K}$ as in Proposition \ref{prop.AoPS333199.gen}. Then, for any
$n\in\mathbb{N}$ and $m\in\mathbb{N}$, we have%
\begin{align}
Y_{m,n}  &  =\sum_{k=0}^{n}\underbrace{y^{k}}_{\substack{=\left(  -1\right)
^{k}\\\text{(since }y=-1\text{)}}}\dbinom{n}{k}\left(  \underbrace{x^{n-k}%
}_{\substack{=X^{n-k}\\\text{(since }x=X\text{)}}}+\underbrace{y}%
_{=-1}\right)  ^{m}\nonumber\\
&  =\sum_{k=0}^{n}\left(  -1\right)  ^{k}\dbinom{n}{k}\left(
\underbrace{X^{n-k}+\left(  -1\right)  }_{=X^{n-k}-1}\right)  ^{m}\nonumber\\
&  =\sum_{k=0}^{n}\left(  -1\right)  ^{k}\dbinom{n}{k}\left(  X^{n-k}%
-1\right)  ^{m}=Z_{m,n}. \label{sol.AoPS333199.1}%
\end{align}

Now, fix $n\in\mathbb{N}$ and $m\in\mathbb{N}$. Applying
(\ref{sol.AoPS333199.1}) to $m$ and $n$ instead of $n$ and $m$, we obtain
$Y_{n,m}=Z_{n,m}$. But Proposition \ref{prop.AoPS333199.gen} shows that
$Y_{m,n}=Y_{n,m}$. Comparing this with (\ref{sol.AoPS333199.1}), we obtain
$Y_{n,m}=Z_{m,n}$. Comparing this with $Y_{n,m}=Z_{n,m}$, we obtain
$Z_{m,n}=Z_{n,m}$. This solves Exercise \ref{exe.AoPS333199}.
\end{proof}

\begin{remark}
Two solutions to Exercise \ref{exe.AoPS333199} are sketched in
\url{http://www.artofproblemsolving.com/community/c6h333199p1782800} . One is
essentially the solution given above (except in lesser generality); the other
is combinatorial.
\end{remark}

\subsection{Solution to Exercise \ref{exe.AoPS262752}}

We shall give two solutions to Exercise \ref{exe.AoPS262752}. The first
solution follows the Hint given in the exercise, and illustrates both the use
of Lemma \ref{lem.polyeq} \textbf{(b)} and of \textquotedblleft generating
functions\textquotedblright\ (the strategy of proving identities by
identifying both sides as coefficients in polynomials or power series). The
second solution is of a more classical nature, using no new methods but a
tricky application of Theorem \ref{thm.vandermonde.XY}.

\subsubsection{First solution}

The crux of the first solution is the proof of the following lemma (which
appears in \cite[(5.55)]{GKP}):

\begin{lemma}
\label{lem.AoPS262752.r}Let $n\in\mathbb{N}$ and $x\in\mathbb{N}$. Then,%
\[
\sum_{k=0}^{n}\left(  -1\right)  ^{k}\dbinom{x}{k}\dbinom{x}{n-k}=%
\begin{cases}
\left(  -1\right)  ^{n/2}\dbinom{x}{n/2}, & \text{if }n\text{ is even};\\
0, & \text{if }n\text{ is odd}%
\end{cases}
.
\]

\end{lemma}

Thus, Lemma \ref{lem.AoPS262752.r} is obtained from Exercise
\ref{exe.AoPS262752} by substituting a nonnegative integer $x$ for the
indeterminate $X$. It thus is clear that Lemma \ref{lem.AoPS262752.r} follows
from Exercise \ref{exe.AoPS262752}. However, for us, the interest lies in the
opposite implication: We shall derive Exercise \ref{exe.AoPS262752} from Lemma
\ref{lem.AoPS262752.r}. Let us, however, prove Lemma \ref{lem.AoPS262752.r}
first. But before we do this, let us state a version of the binomial theorem:

\begin{proposition}
\label{prop.AoPS262752.binom}Let $x\in\mathbb{N}$. Then,%
\[
\left(  1+X\right)  ^{x}=\sum_{k\in\mathbb{N}}\dbinom{x}{k}X^{k}%
\]
(an equality between polynomials in $\mathbb{Z}\left[  X\right]  $). (The sum
$\sum_{k\in\mathbb{N}}\dbinom{x}{k}X^{k}$ is an infinite sum, but only
finitely many of its addends are nonzero, so it is well-defined.)
\end{proposition}

\begin{proof}
[Proof of Proposition \ref{prop.AoPS262752.binom}.]We have%
\begin{align*}
\sum_{k\in\mathbb{N}}\dbinom{x}{k}X^{k}  &  =\underbrace{\sum_{\substack{k\in
\mathbb{N};\\k\leq x}}}_{=\sum_{k=0}^{x}}\dbinom{x}{k}X^{k}+\sum
_{\substack{k\in\mathbb{N};\\k>x}}\underbrace{\dbinom{x}{k}}%
_{\substack{=0\\\text{(by (\ref{eq.binom.0}) (applied to }x\text{ and
}k\\\text{instead of }m\text{ and }n\text{) (since }x<k\text{ (since
}k>x\text{)))}}}X^{k}\\
&  =\sum_{k=0}^{x}\dbinom{x}{k}X^{k}+\underbrace{\sum_{\substack{k\in
\mathbb{N};\\k>x}}0X^{k}}_{=0}=\sum_{k=0}^{x}\dbinom{x}{k}X^{k}.
\end{align*}
Comparing this with
\begin{align*}
\left(  \underbrace{1+X}_{=X+1}\right)  ^{x}  &  =\left(  X+1\right)
^{x}=\sum_{k=0}^{x}\dbinom{x}{k}X^{k}\underbrace{1^{x-k}}_{=1}%
\ \ \ \ \ \ \ \ \ \ \left(  \text{by the binomial formula}\right) \\
&  =\sum_{k=0}^{x}\dbinom{x}{k}X^{k},
\end{align*}
we obtain $\left(  1+X\right)  ^{x}=\sum_{k\in\mathbb{N}}\dbinom{x}{k}X^{k}$.
This proves Proposition \ref{prop.AoPS262752.binom}.
\end{proof}

\begin{proof}
[Proof of Lemma \ref{lem.AoPS262752.r}.]Proposition
\ref{prop.AoPS262752.binom} yields
\begin{equation}
\left(  1+X\right)  ^{x}=\sum_{k\in\mathbb{N}}\dbinom{x}{k}X^{k}.
\label{pf.lem.AoPS262752.r.1}%
\end{equation}
Substituting $-X$ for $X$ in this equality, we obtain%
\[
\left(  1+\left(  -X\right)  \right)  ^{x}=\sum_{k\in\mathbb{N}}\dbinom{x}%
{k}\underbrace{\left(  -X\right)  ^{k}}_{=\left(  -1\right)  ^{k}X^{k}}%
=\sum_{k\in\mathbb{N}}\dbinom{x}{k}\left(  -1\right)  ^{k}X^{k}=\sum
_{k\in\mathbb{N}}\left(  -1\right)  ^{k}\dbinom{x}{k}X^{k}.
\]
Since $1+\left(  -X\right)  =1-X$, this rewrites as%
\[
\left(  1-X\right)  ^{x}=\sum_{k\in\mathbb{N}}\left(  -1\right)  ^{k}%
\dbinom{x}{k}X^{k}.
\]
Multiplying this equality with (\ref{pf.lem.AoPS262752.r.1}), we obtain%
\begin{align*}
\left(  1-X\right)  ^{x}\left(  1+X\right)  ^{x}  &  =\left(  \sum
_{k\in\mathbb{N}}\left(  -1\right)  ^{k}\dbinom{x}{k}X^{k}\right)  \left(
\sum_{k\in\mathbb{N}}\dbinom{x}{k}X^{k}\right) \\
&  =\sum_{k\in\mathbb{N}}\left(  \sum_{i=0}^{k}\left(  -1\right)  ^{i}%
\dbinom{x}{i}\dbinom{x}{k-i}\right)  X^{k}%
\end{align*}
(according to the definition of the product of two polynomials). Hence,
\begin{align}
&  \left(  \text{the coefficient of }X^{n}\text{ in }\left(  1-X\right)
^{x}\left(  1+X\right)  ^{x}\right) \nonumber\\
&  =\sum_{i=0}^{n}\left(  -1\right)  ^{i}\dbinom{x}{i}\dbinom{x}{n-i}%
=\sum_{k=0}^{n}\left(  -1\right)  ^{k}\dbinom{x}{k}\dbinom{x}{n-k}
\label{pf.lem.AoPS262752.r.6}%
\end{align}
(here, we have renamed the summation index $i$ as $k$).

On the other hand,%
\begin{align*}
\left(  1-X\right)  ^{x}\left(  1+X\right)  ^{x}  &  =\left(
\underbrace{\left(  1-X\right)  \left(  1+X\right)  }_{=1-X^{2}=1+\left(
-X^{2}\right)  }\right)  ^{x}=\left(  1+\left(  -X^{2}\right)  \right)
^{x}=\sum_{k\in\mathbb{N}}\dbinom{x}{k}\underbrace{\left(  -X^{2}\right)
^{k}}_{=\left(  -1\right)  ^{k}\left(  X^{2}\right)  ^{k}}\\
&  \ \ \ \ \ \ \ \ \ \ \left(  \text{this follows from substituting }%
-X^{2}\text{ for }X\text{ in (\ref{pf.lem.AoPS262752.r.1})}\right) \\
&  =\sum_{k\in\mathbb{N}}\underbrace{\dbinom{x}{k}\left(  -1\right)  ^{k}%
}_{=\left(  -1\right)  ^{k}\dbinom{x}{k}}\underbrace{\left(  X^{2}\right)
^{k}}_{=X^{2k}}=\sum_{k\in\mathbb{N}}\left(  -1\right)  ^{k}\dbinom{x}%
{k}X^{2k}.
\end{align*}
Hence,%
\[
\left(  \text{the coefficient of }X^{n}\text{ in }\left(  1-X\right)
^{x}\left(  1+X\right)  ^{x}\right)  =%
\begin{cases}
\left(  -1\right)  ^{n/2}\dbinom{x}{n/2}, & \text{if }n\text{ is even};\\
0, & \text{if }n\text{ is odd}%
\end{cases}
.
\]
Comparing this with (\ref{pf.lem.AoPS262752.r.6}), we obtain%
\[
\sum_{k=0}^{n}\left(  -1\right)  ^{k}\dbinom{x}{k}\dbinom{x}{n-k}=%
\begin{cases}
\left(  -1\right)  ^{n/2}\dbinom{x}{n/2}, & \text{if }n\text{ is even};\\
0, & \text{if }n\text{ is odd}%
\end{cases}
.
\]
This proves Lemma \ref{lem.AoPS262752.r}.
\end{proof}

\begin{proof}
[First solution to Exercise \ref{exe.AoPS262752}.]Define two polynomials $P$
and $Q$ (with rational coefficients) by%
\begin{equation}
P=\sum_{k=0}^{n}\left(  -1\right)  ^{k}\dbinom{X}{k}\dbinom{X}{n-k}
\label{sol.AoPS262752.sol1.P}%
\end{equation}
and%
\begin{equation}
Q=%
\begin{cases}
\left(  -1\right)  ^{n/2}\dbinom{X}{n/2}, & \text{if }n\text{ is even};\\
0, & \text{if }n\text{ is odd}%
\end{cases}
. \label{sol.AoPS262752.sol1.Q}%
\end{equation}
For every $x\in\mathbb{N}$, we have%
\begin{align*}
P\left(  x\right)   &  =\sum_{k=0}^{n}\left(  -1\right)  ^{k}\dbinom{x}%
{k}\dbinom{x}{n-k}\ \ \ \ \ \ \ \ \ \ \left(  \text{by the definition of
}P\right) \\
&  =%
\begin{cases}
\left(  -1\right)  ^{n/2}\dbinom{x}{n/2}, & \text{if }n\text{ is even};\\
0, & \text{if }n\text{ is odd}%
\end{cases}
\ \ \ \ \ \ \ \ \ \ \left(  \text{by Lemma \ref{lem.AoPS262752.r}}\right) \\
&  =Q\left(  x\right)  \ \ \ \ \ \ \ \ \ \ \left(  \text{by the definition of
}Q\right)  .
\end{align*}
Hence, Lemma \ref{lem.polyeq} \textbf{(b)} shows that $P=Q$. In light of
(\ref{sol.AoPS262752.sol1.P}) and (\ref{sol.AoPS262752.sol1.Q}), this rewrites
as
\[
\sum_{k=0}^{n}\left(  -1\right)  ^{k}\dbinom{X}{k}\dbinom{X}{n-k}=%
\begin{cases}
\left(  -1\right)  ^{n/2}\dbinom{X}{n/2}, & \text{if }n\text{ is even};\\
0, & \text{if }n\text{ is odd}%
\end{cases}
.
\]
This solves Exercise \ref{exe.AoPS262752}.
\end{proof}

\subsubsection{Second solution}

Now let us prepare for the second solution to Exercise \ref{exe.AoPS262752}.
We shall use the Iverson bracket notation introduced in Definition
\ref{def.iverson}.

We notice that every $n\in\mathbb{N}$ satisfies%
\begin{equation}
\sum_{k=0}^{n}\left(  -1\right)  ^{k}\dbinom{n}{k}=\left[  n=0\right]
\label{sol.AoPS262752.n=0}%
\end{equation}
(by Proposition \ref{prop.binom.1-1}, applied to $m=n$).

Next, we state a simple fact:

\begin{lemma}
\label{lem.sol.AoPS262752.1}Let $m\in\mathbb{N}$ and $i\in\mathbb{N}$. Then,%
\[
\sum_{k=0}^{m}\left(  -1\right)  ^{k}\dbinom{k}{i}\dbinom{m}{k}=\left(
-1\right)  ^{i}\left[  m=i\right]  .
\]

\end{lemma}

\begin{vershort}
\begin{proof}
[Proof of Lemma \ref{lem.sol.AoPS262752.1}.]If $m<i$, then Lemma
\ref{lem.sol.AoPS262752.1} holds\footnote{\textit{Proof.} Assume that $m<i$.
Then, $m\neq i$. Thus, $\left[  m=i\right]  =0$.
\par
But every $k\in\left\{  0,1,\ldots,m\right\}  $ satisfies $k\leq m<i$. Hence,
every $k\in\left\{  0,1,\ldots,m\right\}  $ satisfies $\dbinom{k}{i}=0$ (by
(\ref{eq.binom.0}), applied to $k$ and $i$ instead of $m$ and $n$). Now,
\[
\sum_{k=0}^{m}\left(  -1\right)  ^{k}\underbrace{\dbinom{k}{i}}_{=0}\dbinom
{m}{k}=\sum_{k=0}^{m}\left(  -1\right)  ^{k}0\dbinom{m}{k}=0=\left[
m=i\right]
\]
(since $\left[  m=i\right]  =0$). In other words, Lemma
\ref{lem.sol.AoPS262752.1} holds, qed.}. Hence, for the rest of this proof of
Lemma \ref{lem.sol.AoPS262752.1}, we can WLOG assume that we don't have $m<i$.
Assume this.

We have $m\geq i$ (since we don't have $m<i$). Hence, $m\geq i\geq0$, so that%
\begin{align}
&  \sum_{k=0}^{m}\left(  -1\right)  ^{k}\underbrace{\dbinom{k}{i}\dbinom{m}%
{k}}_{=\dbinom{m}{k}\dbinom{k}{i}}\nonumber\\
&  =\sum_{k=0}^{m}\left(  -1\right)  ^{k}\dbinom{m}{k}\dbinom{k}{i}\nonumber\\
&  =\sum_{k=0}^{i-1}\left(  -1\right)  ^{k}\dbinom{m}{k}\underbrace{\dbinom
{k}{i}}_{\substack{=0\\\text{(by (\ref{eq.binom.0}) (applied to }k\text{ and
}i\\\text{instead of }m\text{ and }n\text{) (since }k<i\text{))}}}+\sum
_{k=i}^{m}\left(  -1\right)  ^{k}\underbrace{\dbinom{m}{k}\dbinom{k}{i}%
}_{\substack{=\dbinom{m}{i}\dbinom{m-i}{k-i}\\\text{(by
(\ref{eq.binom.trinom-rev.m}) (applied to }k\text{ and }i\\\text{instead of
}i\text{ and }a\text{) (since }k\geq i\text{))}}}\nonumber\\
&  =\underbrace{\sum_{k=0}^{i-1}\left(  -1\right)  ^{k}\dbinom{m}{k}0}%
_{=0}+\sum_{k=i}^{m}\left(  -1\right)  ^{k}\dbinom{m}{i}\dbinom{m-i}{k-i}%
=\sum_{k=i}^{m}\left(  -1\right)  ^{k}\dbinom{m}{i}\dbinom{m-i}{k-i}%
\nonumber\\
&  =\sum_{k=0}^{m-i}\underbrace{\left(  -1\right)  ^{k+i}}_{=\left(
-1\right)  ^{k}\left(  -1\right)  ^{i}}\dbinom{m}{i}\dbinom{m-i}{k}\nonumber\\
&  \ \ \ \ \ \ \ \ \ \ \left(  \text{here, we have substituted }k+i\text{ for
}k\text{ in the sum}\right) \nonumber\\
&  =\sum_{k=0}^{m-i}\left(  -1\right)  ^{k}\left(  -1\right)  ^{i}\dbinom
{m}{i}\dbinom{m-i}{k}=\left(  -1\right)  ^{i}\dbinom{m}{i}\underbrace{\sum
_{k=0}^{m-i}\left(  -1\right)  ^{k}\dbinom{m-i}{k}}_{\substack{=\left[
m-i=0\right]  \\\text{(by (\ref{sol.AoPS262752.n=0}) (applied to
}n=m-i\text{))}}}\nonumber\\
&  =\left(  -1\right)  ^{i}\dbinom{m}{i}\left[  \underbrace{m-i=0}%
_{\substack{\text{this is equivalent to}\\m=i}}\right]  =\left(  -1\right)
^{i}\dbinom{m}{i}\left[  m=i\right]  .
\label{pf.lem.sol.AoPS262752.1.short.main}%
\end{align}
But it is easy to see that $\dbinom{m}{i}\left[  m=i\right]  =\left[
m=i\right]  $\ \ \ \ \footnote{\textit{Proof.} We have $\underbrace{\dbinom
{i}{i}}_{=1}\left[  i=i\right]  =\left[  i=i\right]  $. In other words, the
equality $\dbinom{m}{i}\left[  m=i\right]  =\left[  m=i\right]  $ holds in the
case when $m=i$. Therefore, in order to prove this equality, we only need to
consider the case when $m\neq i$. So assume that $m\neq i$. Then, $\left[
m=i\right]  =0$, and thus $\dbinom{m}{i}\underbrace{\left[  m=i\right]  }%
_{=0}=0=\left[  m=i\right]  $, qed.}. Hence,
(\ref{pf.lem.sol.AoPS262752.1.short.main}) becomes%
\[
\sum_{k=0}^{m}\left(  -1\right)  ^{k}\dbinom{k}{i}\dbinom{m}{k}=\left(
-1\right)  ^{i}\underbrace{\dbinom{m}{i}\left[  m=i\right]  }_{=\left[
m=i\right]  }=\left(  -1\right)  ^{i}\left[  m=i\right]  .
\]
This proves Lemma \ref{lem.sol.AoPS262752.1}.
\end{proof}
\end{vershort}

\begin{verlong}
\begin{proof}
[Proof of Lemma \ref{lem.sol.AoPS262752.1}.]If $m<i$, then Lemma
\ref{lem.sol.AoPS262752.1} holds\footnote{\textit{Proof.} Assume that $m<i$.
Then, $m\neq i$. Hence, the statement $m=i$ is false. Thus, $\left[
m=i\right]  =0$.
\par
But every $k\in\left\{  0,1,\ldots,m\right\}  $ satisfies $k\leq m<i$. Hence,
every $k\in\left\{  0,1,\ldots,m\right\}  $ satisfies
\begin{equation}
\dbinom{k}{i}=0 \label{pf.lem.sol.AoPS262752.1.fn1.1}%
\end{equation}
(by (\ref{eq.binom.0}), applied to $k$ and $i$ instead of $m$ and $n$). Now,
\[
\sum_{k=0}^{m}\left(  -1\right)  ^{k}\underbrace{\dbinom{k}{i}}%
_{\substack{=0\\\text{(by (\ref{pf.lem.sol.AoPS262752.1.fn1.1}))}}}\dbinom
{m}{k}=\sum_{k=0}^{m}\left(  -1\right)  ^{k}0\dbinom{m}{k}=0=\left[
m=i\right]
\]
(since $\left[  m=i\right]  =0$). In other words, Lemma
\ref{lem.sol.AoPS262752.1} holds, qed.}. Hence, for the rest of this proof of
Lemma \ref{lem.sol.AoPS262752.1}, we can WLOG assume that we don't have $m<i$.
Assume this.

We have $m\geq i$ (since we don't have $m<i$). Hence, $m\geq i\geq0$, so that%
\begin{align}
&  \sum_{k=0}^{m}\left(  -1\right)  ^{k}\underbrace{\dbinom{k}{i}\dbinom{m}%
{k}}_{=\dbinom{m}{k}\dbinom{k}{i}}\nonumber\\
&  =\sum_{k=0}^{m}\left(  -1\right)  ^{k}\dbinom{m}{k}\dbinom{k}{i}\nonumber\\
&  =\sum_{k=0}^{i-1}\left(  -1\right)  ^{k}\dbinom{m}{k}\underbrace{\dbinom
{k}{i}}_{\substack{=0\\\text{(by (\ref{eq.binom.0}) (applied to }k\text{ and
}i\\\text{instead of }m\text{ and }n\text{) (since }k<i\text{))}}}+\sum
_{k=i}^{m}\left(  -1\right)  ^{k}\underbrace{\dbinom{m}{k}\dbinom{k}{i}%
}_{\substack{=\dbinom{m}{i}\dbinom{m-i}{k-i}\\\text{(by
(\ref{eq.binom.trinom-rev.m}) (applied to }k\text{ and }i\\\text{instead of
}i\text{ and }a\text{) (since }k\geq i\text{))}}}\nonumber\\
&  =\underbrace{\sum_{k=0}^{i-1}\left(  -1\right)  ^{k}\dbinom{m}{k}0}%
_{=0}+\sum_{k=i}^{m}\left(  -1\right)  ^{k}\dbinom{m}{i}\dbinom{m-i}{k-i}%
=\sum_{k=i}^{m}\left(  -1\right)  ^{k}\dbinom{m}{i}\dbinom{m-i}{k-i}%
\nonumber\\
&  =\sum_{k=0}^{m-i}\underbrace{\left(  -1\right)  ^{k+i}}_{=\left(
-1\right)  ^{k}\left(  -1\right)  ^{i}}\dbinom{m}{i}\underbrace{\dbinom
{m-i}{k+i-i}}_{\substack{=\dbinom{m-i}{k}\\\text{(since }k+i-i=k\text{)}%
}}\nonumber\\
&  \ \ \ \ \ \ \ \ \ \ \left(  \text{here, we have substituted }k+i\text{ for
}k\text{ in the sum}\right) \nonumber\\
&  =\sum_{k=0}^{m-i}\left(  -1\right)  ^{k}\left(  -1\right)  ^{i}\dbinom
{m}{i}\dbinom{m-i}{k}=\left(  -1\right)  ^{i}\dbinom{m}{i}\underbrace{\sum
_{k=0}^{m-i}\left(  -1\right)  ^{k}\dbinom{m-i}{k}}_{\substack{=\left[
m-i=0\right]  \\\text{(by (\ref{sol.AoPS262752.n=0}) (applied to
}n=m-i\text{))}}}\nonumber\\
&  =\left(  -1\right)  ^{i}\dbinom{m}{i}\left[  m-i=0\right]  .
\label{pf.lem.sol.AoPS262752.1.main}%
\end{align}

If $m=i$, then Lemma \ref{lem.sol.AoPS262752.1} holds\footnote{\textit{Proof.}
Assume that $m=i$. Hence, $m-i=0$, so that $\left[  m-i=0\right]  =1$. Also,
$\left[  m=i\right]  =1$ (since $m=i$). Also, from $m=i$, we obtain
$\dbinom{m}{i}=\dbinom{i}{i}=1$ (by (\ref{eq.binom.mm}), applied to $m=i$).
Hence, (\ref{pf.lem.sol.AoPS262752.1.main}) becomes%
\[
\sum_{k=0}^{m}\left(  -1\right)  ^{k}\dbinom{k}{i}\dbinom{m}{k}=\left(
-1\right)  ^{i}\underbrace{\dbinom{m}{i}}_{=1}\underbrace{\left[
m-i=0\right]  }_{=1=\left[  m=i\right]  }=\left(  -1\right)  ^{i}\left[
m=i\right]  .
\]
In other words, Lemma \ref{lem.sol.AoPS262752.1} holds, qed.}. Hence, for the
rest of this proof of Lemma \ref{lem.sol.AoPS262752.1}, we can WLOG assume
that we don't have $m=i$. Assume this.

We have $m\neq i$ (since we don't have $m=i$). Hence, $\left[  m=i\right]
=0$, so that $\left(  -1\right)  ^{i}\underbrace{\left[  m=i\right]  }_{=0}%
=0$. Also, $m-i\neq0$ (since $m\neq i$) and thus $\left[  m-i=0\right]  =0$.
Now, (\ref{pf.lem.sol.AoPS262752.1.main}) becomes%
\[
\sum_{k=0}^{m}\left(  -1\right)  ^{k}\dbinom{k}{i}\dbinom{m}{k}=\left(
-1\right)  ^{i}\dbinom{m}{i}\underbrace{\left[  m-i=0\right]  }_{=0}=0=\left(
-1\right)  ^{i}\left[  m=i\right]
\]
(since $\left(  -1\right)  ^{i}\left[  m=i\right]  =0$). This proves Lemma
\ref{lem.sol.AoPS262752.1}.
\end{proof}
\end{verlong}

Here comes one more simple lemma:

\begin{lemma}
\label{lem.sol.AoPS262752.2}Let $n\in\mathbb{N}$. Let $a_{0},a_{1}%
,\ldots,a_{n}$ be $n+1$ polynomials in the indeterminate $X$ with rational
coefficients (that is, $n+1$ elements of $\mathbb{Q}\left[  X\right]  $).
Then,%
\[
\sum_{i=0}^{n}a_{i}\left[  n-i=i\right]  =%
\begin{cases}
a_{n/2}, & \text{if }n\text{ is even};\\
0, & \text{if }n\text{ is odd}%
\end{cases}
.
\]

\end{lemma}

\begin{proof}
[Proof of Lemma \ref{lem.sol.AoPS262752.2}.]We have $n\in\mathbb{N}$, so that
$n\geq0$. For every $i\in\left\{  0,1,\ldots,n\right\}  $, we have%
\begin{equation}
\left[  \underbrace{n-i=i}_{\text{this is equivalent to }n=2i}\right]
=\left[  \underbrace{n=2i}_{\text{this is equivalent to }i=n/2}\right]
=\left[  i=n/2\right]  . \label{pf.lem.sol.AoPS262752.2.1}%
\end{equation}
Thus,%
\begin{align}
&  \underbrace{\sum_{i=0}^{n}}_{=\sum_{i\in\left\{  0,1,\ldots,n\right\}  }%
}a_{i}\underbrace{\left[  n-i=i\right]  }_{\substack{=\left[  i=n/2\right]
\\\text{(by (\ref{pf.lem.sol.AoPS262752.2.1}))}}}\nonumber\\
&  =\sum_{i\in\left\{  0,1,\ldots,n\right\}  }a_{i}\left[  i=n/2\right]
=\sum_{\substack{i\in\left\{  0,1,\ldots,n\right\}  ;\\i=n/2}}a_{i}%
\underbrace{\left[  i=n/2\right]  }_{\substack{=1\\\text{(since }%
i=n/2\text{)}}}+\sum_{\substack{i\in\left\{  0,1,\ldots,n\right\}  ;\\i\neq
n/2}}a_{i}\underbrace{\left[  i=n/2\right]  }_{\substack{=0\\\text{(since
}i\neq n/2\text{)}}}\nonumber\\
&  \ \ \ \ \ \ \ \ \ \ \ \ \ \ \ \ \ \ \ \ \left(
\begin{array}
[c]{c}%
\text{since every }i\in\left\{  0,1,\ldots,n\right\}  \text{ satisfies}\\
\text{either }i=n/2\text{ or }i\neq n/2\text{ (but not both)}%
\end{array}
\right) \nonumber\\
&  =\sum_{\substack{i\in\left\{  0,1,\ldots,n\right\}  ;\\i=n/2}%
}a_{i}+\underbrace{\sum_{\substack{i\in\left\{  0,1,\ldots,n\right\}  ;\\i\neq
n/2}}a_{i}0}_{=0}=\sum_{\substack{i\in\left\{  0,1,\ldots,n\right\}
;\\i=n/2}}a_{i}. \label{pf.lem.sol.AoPS262752.2.2}%
\end{align}

We must be in one of the following two cases:

\textit{Case 1:} The number $n$ is even.

\textit{Case 2:} The number $n$ is odd.

Let us first consider Case 1. In this case, the number $n$ is even. Hence,
\begin{equation}%
\begin{cases}
a_{n/2}, & \text{if }n\text{ is even};\\
0, & \text{if }n\text{ is odd}%
\end{cases}
=a_{n/2}. \label{pf.lem.sol.AoPS262752.2.c1.1}%
\end{equation}

On the other hand, $n/2\in\mathbb{Z}$ (since $n$ is even). Combined with
$n/2\geq0$ (since $n\geq0$) and $n/2\leq n$ (for the same reason), this shows
that $n/2\in\left\{  0,1,\ldots,n\right\}  $. Now,
(\ref{pf.lem.sol.AoPS262752.2.2}) becomes%
\begin{align*}
\sum_{i=0}^{n}a_{i}\left[  n-i=i\right]   &  =\sum_{\substack{i\in\left\{
0,1,\ldots,n\right\}  ;\\i=n/2}}a_{i}=a_{n/2}\ \ \ \ \ \ \ \ \ \ \left(
\text{since }n/2\in\left\{  0,1,\ldots,n\right\}  \right) \\
&  =%
\begin{cases}
a_{n/2}, & \text{if }n\text{ is even};\\
0, & \text{if }n\text{ is odd}%
\end{cases}
\ \ \ \ \ \ \ \ \ \ \left(  \text{by (\ref{pf.lem.sol.AoPS262752.2.c1.1}%
)}\right)  .
\end{align*}
Thus, Lemma \ref{lem.sol.AoPS262752.2} is proven in Case 1.

Let us now consider Case 2. In this case, the number $n$ is odd. Hence,
\begin{equation}%
\begin{cases}
a_{n/2}, & \text{if }n\text{ is even};\\
0, & \text{if }n\text{ is odd}%
\end{cases}
=0. \label{pf.lem.sol.AoPS262752.2.c2.1}%
\end{equation}

On the other hand, $n/2\notin\mathbb{Z}$ (since $n$ is odd). Hence,
$n/2\notin\left\{  0,1,\ldots,n\right\}  $ (since $\left\{  0,1,\ldots
,n\right\}  \subseteq\mathbb{Z}$). Now, (\ref{pf.lem.sol.AoPS262752.2.2})
becomes%
\begin{align*}
\sum_{i=0}^{n}a_{i}\left[  n-i=i\right]   &  =\sum_{\substack{i\in\left\{
0,1,\ldots,n\right\}  ;\\i=n/2}}a_{i}=\left(  \text{empty sum}\right)
\ \ \ \ \ \ \ \ \ \ \left(  \text{since }n/2\notin\left\{  0,1,\ldots
,n\right\}  \right) \\
&  =0=%
\begin{cases}
a_{n/2}, & \text{if }n\text{ is even};\\
0, & \text{if }n\text{ is odd}%
\end{cases}
\ \ \ \ \ \ \ \ \ \ \left(  \text{by (\ref{pf.lem.sol.AoPS262752.2.c2.1}%
)}\right)  .
\end{align*}
Thus, Lemma \ref{lem.sol.AoPS262752.2} is proven in Case 2.

We have now proved Lemma \ref{lem.sol.AoPS262752.2} in both Cases 1 and 2.
Since these two Cases cover all possibilities, this shows that Lemma
\ref{lem.sol.AoPS262752.2} always holds.
\end{proof}

Now, we are ready to solve Exercise \ref{exe.AoPS262752} again:

\begin{proof}
[Second solution to Exercise \ref{exe.AoPS262752}.]Let $g\in\left\{
0,1,\ldots,n\right\}  $ be arbitrary. Then, $n-g\in\left\{  0,1,\ldots
,n\right\}  \subseteq\mathbb{N}$. Hence, Theorem \ref{thm.vandermonde.XY}
(applied to $n-g$ instead of $n$) yields%
\[
\dbinom{X+Y}{n-g}=\sum_{k=0}^{n-g}\dbinom{X}{k}\dbinom{Y}{n-g-k}%
\]
(an equality between two polynomials in $X$ and $Y$). Substituting $g$ and
$X-g$ for $X$ and $Y$ in this equality, we obtain%
\[
\dbinom{g+\left(  X-g\right)  }{n-g}=\sum_{k=0}^{n-g}\dbinom{g}{k}\dbinom
{X-g}{n-g-k}=\sum_{i=0}^{n-g}\dbinom{g}{i}\dbinom{X-g}{n-g-i}%
\]
(here, we have renamed the summation index $k$ as $i$). Since $g+\left(
X-g\right)  =X$, this rewrites as
\begin{equation}
\dbinom{X}{n-g}=\sum_{i=0}^{n-g}\dbinom{g}{i}\dbinom{X-g}{n-g-i}.
\label{sol.AoPS262752.sol2.1}%
\end{equation}

Now, let us forget that we fixed $g$. We thus have shown that
(\ref{sol.AoPS262752.sol2.1}) holds for every $g\in\left\{  0,1,\ldots
,n\right\}  $.

On the other hand, for every $k\in\mathbb{N}$ and $i\in\mathbb{N}$ satisfying
$k+i\leq n$, we have%
\begin{equation}
\dbinom{X}{k}\dbinom{X-k}{n-k-i}=\dbinom{X}{n-i}\dbinom{n-i}{k}
\label{sol.AoPS262752.sol2.3}%
\end{equation}
\footnote{\textit{Proof of (\ref{sol.AoPS262752.sol2.3}):} Let $k\in
\mathbb{N}$ and $i\in\mathbb{N}$ be such that $k+i\leq n$. From $k+i\leq n$,
we obtain $n-i\geq k$. Thus, $n-i\geq k\geq0$, so that $n-i\in\mathbb{N}$.
Hence, Proposition \ref{prop.binom.Xes} \textbf{(f)} (applied to $n-i$ and $k$
instead of $i$ and $a$) shows that
\[
\dbinom{X}{n-i}\dbinom{n-i}{k}=\dbinom{X}{k}\dbinom{X-k}{\left(  n-i\right)
-k}=\dbinom{X}{k}\dbinom{X-k}{n-k-i}%
\]
(since $\left(  n-i\right)  -k=n-k-i$). This proves
(\ref{sol.AoPS262752.sol2.3}).}.

Now,%
\begin{align}
&  \sum_{k=0}^{n}\left(  -1\right)  ^{k}\dbinom{X}{k}\underbrace{\dbinom
{X}{n-k}}_{\substack{=\sum_{i=0}^{n-k}\dbinom{k}{i}\dbinom{X-k}{n-k-i}%
\\\text{(by (\ref{sol.AoPS262752.sol2.1}) (applied to }g=k\text{))}%
}}\nonumber\\
&  =\sum_{k=0}^{n}\left(  -1\right)  ^{k}\dbinom{X}{k}\left(  \sum_{i=0}%
^{n-k}\dbinom{k}{i}\dbinom{X-k}{n-k-i}\right) \nonumber\\
&  =\underbrace{\sum_{k=0}^{n}\sum_{i=0}^{n-k}}_{\substack{=\sum
_{\substack{\left(  k,i\right)  \in\mathbb{N}^{2};\\k\leq n;\ i\leq n-k}%
}=\sum_{\substack{\left(  k,i\right)  \in\mathbb{N}^{2};\\k\leq n;\ k+i\leq
n}}\\\text{(since the condition }i\leq n-k\text{ is equivalent to }k+i\leq
n\text{)}}}\left(  -1\right)  ^{k}\underbrace{\dbinom{X}{k}\dbinom{k}{i}%
}_{=\dbinom{k}{i}\dbinom{X}{k}}\dbinom{X-k}{n-k-i}\nonumber\\
&  =\underbrace{\sum_{\substack{\left(  k,i\right)  \in\mathbb{N}^{2};\\k\leq
n;\ k+i\leq n}}}_{\substack{=\sum_{\substack{\left(  k,i\right)  \in
\mathbb{N}^{2};\\k+i\leq n}}\\\text{(since the condition }\left(  k\leq
n\text{ and }k+i\leq n\right)  \text{ is equivalent to }k+i\leq
n\\\text{(because the condition }k\leq n\text{ follows from }k+i\leq
n\text{))}}}\left(  -1\right)  ^{k}\dbinom{k}{i}\underbrace{\dbinom{X}%
{k}\dbinom{X-k}{n-k-i}}_{\substack{=\dbinom{X}{n-i}\dbinom{n-i}{k}\\\text{(by
(\ref{sol.AoPS262752.sol2.3}))}}}\nonumber\\
&  =\underbrace{\sum_{\substack{\left(  k,i\right)  \in\mathbb{N}%
^{2};\\k+i\leq n}}}_{\substack{=\sum_{\substack{\left(  k,i\right)
\in\mathbb{N}^{2};\\i\leq n;\ k+i\leq n}}\\\text{(since the condition }k+i\leq
n\text{ is equivalent to }\left(  i\leq n\text{ and }k+i\leq n\right)
\\\text{(because the condition }i\leq n\text{ follows from }k+i\leq
n\text{))}}}\left(  -1\right)  ^{k}\dbinom{k}{i}\dbinom{X}{n-i}\dbinom{n-i}%
{k}\nonumber
\end{align}%
\begin{align*}
&  =\underbrace{\sum_{\substack{\left(  k,i\right)  \in\mathbb{N}^{2};\\i\leq
n;\ k+i\leq n}}}_{\substack{=\sum_{\substack{\left(  i,k\right)  \in
\mathbb{N}^{2};\\i\leq n;\ k+i\leq n}}=\sum_{\substack{\left(  i,k\right)
\in\mathbb{N}^{2};\\i\leq n;\ k\leq n-i}}\\\text{(since the condition }k+i\leq
n\text{ is equivalent to }k\leq n-i\text{)}}}\left(  -1\right)  ^{k}\dbinom
{k}{i}\dbinom{X}{n-i}\dbinom{n-i}{k}\\
&  =\underbrace{\sum_{\substack{\left(  i,k\right)  \in\mathbb{N}^{2};\\i\leq
n;\ k\leq n-i}}}_{=\sum_{i=0}^{n}\sum_{k=0}^{n-i}}\left(  -1\right)
^{k}\dbinom{k}{i}\dbinom{X}{n-i}\dbinom{n-i}{k}\\
&  =\sum_{i=0}^{n}\sum_{k=0}^{n-i}\left(  -1\right)  ^{k}\dbinom{k}{i}%
\dbinom{X}{n-i}\dbinom{n-i}{k}=\sum_{i=0}^{n}\dbinom{X}{n-i}\underbrace{\sum
_{k=0}^{n-i}\left(  -1\right)  ^{k}\dbinom{k}{i}\dbinom{n-i}{k}}%
_{\substack{=\left(  -1\right)  ^{i}\left[  n-i=i\right]  \\\text{(by Lemma
\ref{lem.sol.AoPS262752.1}}\\\text{(applied to }m=n-i\text{))}}}\\
&  =\sum_{i=0}^{n}\underbrace{\dbinom{X}{n-i}\left(  -1\right)  ^{i}%
}_{=\left(  -1\right)  ^{i}\dbinom{X}{n-i}}\left[  n-i=i\right] \\
&  =\sum_{i=0}^{n}\left(  -1\right)  ^{i}\dbinom{X}{n-i}\left[  n-i=i\right]
=%
\begin{cases}
\left(  -1\right)  ^{n/2}\dbinom{X}{n-n/2}, & \text{if }n\text{ is even};\\
0, & \text{if }n\text{ is odd}%
\end{cases}
\\
&  \ \ \ \ \ \ \ \ \ \ \ \ \ \ \ \ \ \ \ \ \left(  \text{by Lemma
\ref{lem.sol.AoPS262752.2}, applied to }a_{i}=\left(  -1\right)  ^{i}%
\dbinom{X}{n-i}\right) \\
&  =%
\begin{cases}
\left(  -1\right)  ^{n/2}\dbinom{X}{n/2}, & \text{if }n\text{ is even};\\
0, & \text{if }n\text{ is odd}%
\end{cases}
\ \ \ \ \ \ \ \ \ \ \left(  \text{since }n-n/2=n/2\right)  .
\end{align*}
This solves Exercise \ref{exe.AoPS262752}.
\end{proof}

\subsubsection{Addendum}

Let us record a classical result which follows from Exercise
\ref{exe.AoPS262752}:

\begin{corollary}
\label{cor.AoPS262752.X=n}Let $n\in\mathbb{N}$. Then,%
\[
\sum_{k=0}^{n}\left(  -1\right)  ^{k}\dbinom{n}{k}^{2}=%
\begin{cases}
\left(  -1\right)  ^{n/2}\dbinom{n}{n/2}, & \text{if }n\text{ is even};\\
0, & \text{if }n\text{ is odd}%
\end{cases}
.
\]

\end{corollary}

\begin{proof}
[Proof of Corollary \ref{cor.AoPS262752.X=n}.]Exercise \ref{exe.AoPS262752}
shows that%
\[
\sum_{k=0}^{n}\left(  -1\right)  ^{k}\dbinom{X}{k}\dbinom{X}{n-k}=%
\begin{cases}
\left(  -1\right)  ^{n/2}\dbinom{X}{n/2}, & \text{if }n\text{ is even};\\
0, & \text{if }n\text{ is odd}%
\end{cases}
.
\]
Substituting $n$ for $X$ in this equality, we obtain%
\[
\sum_{k=0}^{n}\left(  -1\right)  ^{k}\dbinom{n}{k}\dbinom{n}{n-k}=%
\begin{cases}
\left(  -1\right)  ^{n/2}\dbinom{n}{n/2}, & \text{if }n\text{ is even};\\
0, & \text{if }n\text{ is odd}%
\end{cases}
.
\]
Now,%
\begin{align*}
\sum_{k=0}^{n}\left(  -1\right)  ^{k}\underbrace{\dbinom{n}{k}^{2}}%
_{=\dbinom{n}{k}\dbinom{n}{k}}  &  =\sum_{k=0}^{n}\left(  -1\right)
^{k}\dbinom{n}{k}\underbrace{\dbinom{n}{k}}_{\substack{=\dbinom{n}%
{n-k}\\\text{(by (\ref{eq.binom.symm}) (applied to }n\text{ and }%
k\\\text{instead of }m\text{ and }n\text{))}}}\\
&  =\sum_{k=0}^{n}\left(  -1\right)  ^{k}\dbinom{n}{k}\dbinom{n}{n-k}=%
\begin{cases}
\left(  -1\right)  ^{n/2}\dbinom{n}{n/2}, & \text{if }n\text{ is even};\\
0, & \text{if }n\text{ is odd}%
\end{cases}
.
\end{align*}
This proves Corollary \ref{cor.AoPS262752.X=n}.
\end{proof}

Corollary \ref{cor.AoPS262752.X=n} is a well-known fact. Mike Spivey has found
a combinatorial proof \cite{Spivey12} (which generalizes immediately to a
proof of Lemma \ref{lem.AoPS262752.r}); the corollary also has appeared in
\url{http://www.artofproblemsolving.com/community/c6h262752} .

\subsection{Solution to Exercise \ref{exe.vander-1/2}}

Exercise \ref{exe.vander-1/2} \textbf{(a)} is a known fact. It appears, e.g.,
in \cite[\S 5.3, (5.39)]{GKP}. Combinatorial proofs appear in \cite{Sved84} as
well as in the math.stackexchange discussions
\url{https://math.stackexchange.com/questions/72367} and
\url{https://math.stackexchange.com/a/360780} . Complicated combinatorial
proofs of Exercise \ref{exe.vander-1/2} \textbf{(b)} have been given in
\cite{Spivey12b} and at \url{https://math.stackexchange.com/questions/80649} .
We shall prove both parts of Exercise \ref{exe.vander-1/2} algebraically,
making use of Exercise \ref{exe.bin.-1/2} \textbf{(b)}. (This is how part
\textbf{(a)} of this exercise is proven in \cite[\S 5.3]{GKP}.)

\begin{proof}
[Solution to Exercise \ref{exe.vander-1/2}.]Let $k\in\left\{  0,1,\ldots
,n\right\}  $. Thus, $n-k\in\left\{  0,1,\ldots,n\right\}  \subseteq
\mathbb{N}$. Hence, Exercise \ref{exe.bin.-1/2} \textbf{(b)} (applied to $n-k$
instead of $n$) yields%
\begin{equation}
\dbinom{-1/2}{n-k}=\left(  \dfrac{-1}{4}\right)  ^{n-k}\dbinom{2\left(
n-k\right)  }{n-k}. \label{sol.vander-1/2.factor2}%
\end{equation}
But $k\in\left\{  0,1,\ldots,n\right\}  \subseteq\mathbb{N}$. Hence, Exercise
\ref{exe.bin.-1/2} \textbf{(b)} (applied to $k$ instead of $n$) yields%
\begin{equation}
\dbinom{-1/2}{k}=\left(  \dfrac{-1}{4}\right)  ^{k}\dbinom{2k}{k}.
\label{sol.vander-1/2.factor1}%
\end{equation}
Multiplying the equalities (\ref{sol.vander-1/2.factor1}) and
(\ref{sol.vander-1/2.factor2}), we obtain%
\begin{align}
\dbinom{-1/2}{k}\dbinom{-1/2}{n-k}  &  =\left(  \dfrac{-1}{4}\right)
^{k}\dbinom{2k}{k}\left(  \dfrac{-1}{4}\right)  ^{n-k}\dbinom{2\left(
n-k\right)  }{n-k}\nonumber\\
&  =\underbrace{\left(  \dfrac{-1}{4}\right)  ^{k}\left(  \dfrac{-1}%
{4}\right)  ^{n-k}}_{\substack{=\left(  \dfrac{-1}{4}\right)  ^{k+\left(
n-k\right)  }=\left(  \dfrac{-1}{4}\right)  ^{n}\\\text{(since }k+\left(
n-k\right)  =n\text{)}}}\dbinom{2k}{k}\dbinom{2\left(  n-k\right)  }%
{n-k}\nonumber\\
&  =\left(  \dfrac{-1}{4}\right)  ^{n}\dbinom{2k}{k}\dbinom{2\left(
n-k\right)  }{n-k}. \label{sol.vander-1/2.0}%
\end{align}

Now, forget that we fixed $k$. We thus have proven the equality
(\ref{sol.vander-1/2.0}) for each $k\in\left\{  0,1,\ldots,n\right\}  $.

We also recall that $a^{n}b^{n}=\left(  ab\right)  ^{n}$ for any
$a,b\in\mathbb{Q}$. Applying this to $a=-4$ and $b=\dfrac{-1}{4}$, we obtain%
\[
\left(  -4\right)  ^{n}\left(  \dfrac{-1}{4}\right)  ^{n}=\left(
\underbrace{\left(  -4\right)  \left(  \dfrac{-1}{4}\right)  }_{=1}\right)
^{n}=1^{n}=1.
\]

\textbf{(a)} Theorem \ref{thm.vandermonde.rat} (applied to $x=-1/2$ and
$y=-1/2$) yields%
\[
\dbinom{\left(  -1/2\right)  +\left(  -1/2\right)  }{n}=\sum_{k=0}^{n}%
\dbinom{-1/2}{k}\dbinom{-1/2}{n-k}.
\]
Comparing this with%
\begin{align*}
\dbinom{\left(  -1/2\right)  +\left(  -1/2\right)  }{n}  &  =\dbinom{-1}%
{n}\ \ \ \ \ \ \ \ \ \ \left(  \text{since }\left(  -1/2\right)  +\left(
-1/2\right)  =-1\right) \\
&  =\left(  -1\right)  ^{n}\ \ \ \ \ \ \ \ \ \ \left(  \text{by Corollary
\ref{cor.binom.-1}}\right)  ,
\end{align*}
we obtain%
\begin{align*}
\left(  -1\right)  ^{n}  &  =\sum_{k=0}^{n}\underbrace{\dbinom{-1/2}{k}%
\dbinom{-1/2}{n-k}}_{\substack{=\left(  \dfrac{-1}{4}\right)  ^{n}\dbinom
{2k}{k}\dbinom{2\left(  n-k\right)  }{n-k}\\\text{(by (\ref{sol.vander-1/2.0}%
))}}}=\sum_{k=0}^{n}\left(  \dfrac{-1}{4}\right)  ^{n}\dbinom{2k}{k}%
\dbinom{2\left(  n-k\right)  }{n-k}\\
&  =\left(  \dfrac{-1}{4}\right)  ^{n}\sum_{k=0}^{n}\dbinom{2k}{k}%
\dbinom{2\left(  n-k\right)  }{n-k}.
\end{align*}
Multiplying both sides of this equality by $\left(  -4\right)  ^{n}$, we
obtain%
\[
\left(  -4\right)  ^{n}\left(  -1\right)  ^{n}=\underbrace{\left(  -4\right)
^{n}\left(  \dfrac{-1}{4}\right)  ^{n}}_{=1}\sum_{k=0}^{n}\dbinom{2k}%
{k}\dbinom{2\left(  n-k\right)  }{n-k}=\sum_{k=0}^{n}\dbinom{2k}{k}%
\dbinom{2\left(  n-k\right)  }{n-k}.
\]
Hence,%
\[
\sum_{k=0}^{n}\dbinom{2k}{k}\dbinom{2\left(  n-k\right)  }{n-k}=\left(
-4\right)  ^{n}\left(  -1\right)  ^{n}=\left(  \underbrace{\left(  -4\right)
\left(  -1\right)  }_{=4}\right)  ^{n}=4^{n}.
\]
This solves Exercise \ref{exe.vander-1/2} \textbf{(a)}.

\textbf{(b)} Exercise \ref{exe.AoPS262752} yields%
\[
\sum_{k=0}^{n}\left(  -1\right)  ^{k}\dbinom{X}{k}\dbinom{X}{n-k}=%
\begin{cases}
\left(  -1\right)  ^{n/2}\dbinom{X}{n/2}, & \text{if }n\text{ is even};\\
0, & \text{if }n\text{ is odd}%
\end{cases}
\]
(an identity between polynomials in $\mathbb{Q}\left[  X\right]  $). If we
substitute $-1/2$ for $X$ in this equality, we obtain%
\[
\sum_{k=0}^{n}\left(  -1\right)  ^{k}\dbinom{-1/2}{k}\dbinom{-1/2}{n-k}=%
\begin{cases}
\left(  -1\right)  ^{n/2}\dbinom{-1/2}{n/2}, & \text{if }n\text{ is even};\\
0, & \text{if }n\text{ is odd}%
\end{cases}
.
\]
Hence,%
\begin{align*}%
\begin{cases}
\left(  -1\right)  ^{n/2}\dbinom{-1/2}{n/2}, & \text{if }n\text{ is even};\\
0, & \text{if }n\text{ is odd}%
\end{cases}
&  =\sum_{k=0}^{n}\left(  -1\right)  ^{k}\underbrace{\dbinom{-1/2}{k}%
\dbinom{-1/2}{n-k}}_{\substack{=\left(  \dfrac{-1}{4}\right)  ^{n}\dbinom
{2k}{k}\dbinom{2\left(  n-k\right)  }{n-k}\\\text{(by (\ref{sol.vander-1/2.0}%
))}}}\\
&  =\sum_{k=0}^{n}\left(  -1\right)  ^{k}\left(  \dfrac{-1}{4}\right)
^{n}\dbinom{2k}{k}\dbinom{2\left(  n-k\right)  }{n-k}\\
&  =\left(  \dfrac{-1}{4}\right)  ^{n}\sum_{k=0}^{n}\left(  -1\right)
^{k}\dbinom{2k}{k}\dbinom{2\left(  n-k\right)  }{n-k}.
\end{align*}
Multiplying both sides of this equality by $\left(  -4\right)  ^{n}$, we
obtain%
\begin{align*}
&  \left(  -4\right)  ^{n}%
\begin{cases}
\left(  -1\right)  ^{n/2}\dbinom{-1/2}{n/2}, & \text{if }n\text{ is even};\\
0, & \text{if }n\text{ is odd}%
\end{cases}
\\
&  =\underbrace{\left(  -4\right)  ^{n}\left(  \dfrac{-1}{4}\right)  ^{n}%
}_{=1}\sum_{k=0}^{n}\left(  -1\right)  ^{k}\dbinom{2k}{k}\dbinom{2\left(
n-k\right)  }{n-k}=\sum_{k=0}^{n}\left(  -1\right)  ^{k}\dbinom{2k}{k}%
\dbinom{2\left(  n-k\right)  }{n-k}.
\end{align*}
Hence,%
\begin{align}
&  \sum_{k=0}^{n}\left(  -1\right)  ^{k}\dbinom{2k}{k}\dbinom{2\left(
n-k\right)  }{n-k}\nonumber\\
&  =\left(  -4\right)  ^{n}%
\begin{cases}
\left(  -1\right)  ^{n/2}\dbinom{-1/2}{n/2}, & \text{if }n\text{ is even};\\
0, & \text{if }n\text{ is odd}%
\end{cases}
\nonumber\\
&  =%
\begin{cases}
\left(  -4\right)  ^{n}\cdot\left(  -1\right)  ^{n/2}\dbinom{-1/2}{n/2}, &
\text{if }n\text{ is even};\\
\left(  -4\right)  ^{n}\cdot0, & \text{if }n\text{ is odd}%
\end{cases}
\nonumber\\
&  =%
\begin{cases}
\left(  -4\right)  ^{n}\cdot\left(  -1\right)  ^{n/2}\dbinom{-1/2}{n/2}, &
\text{if }n\text{ is even};\\
0, & \text{if }n\text{ is odd}%
\end{cases}
\label{sol.vander-1/2.b.step2}%
\end{align}
(since $\left(  -4\right)  ^{n}\cdot0=0$).

Now, let us assume that $n$ is even. Thus, $n/2\in\mathbb{N}$. Hence,
$2^{n}=2^{2\left(  n/2\right)  }=\left(  2^{2}\right)  ^{n/2}=4^{n/2}$ (since
$2^{2}=4$). From $n-n/2=n/2$, we obtain $4^{n-n/2}=4^{n/2}=2^{n}$. But recall
that $n/2\in\mathbb{N}$. Hence, Exercise \ref{exe.bin.-1/2} \textbf{(b)}
(applied to $n/2$ instead of $n$) yields%
\[
\dbinom{-1/2}{n/2}=\underbrace{\left(  \dfrac{-1}{4}\right)  ^{n/2}}%
_{=\dfrac{\left(  -1\right)  ^{n/2}}{4^{n/2}}}\ \ \underbrace{\dbinom{2\left(
n/2\right)  }{n/2}}_{\substack{=\dbinom{n}{n/2}\\\text{(since }2\left(
n/2\right)  =n\text{)}}}=\dfrac{\left(  -1\right)  ^{n/2}}{4^{n/2}}\dbinom
{n}{n/2}.
\]
Hence,%
\begin{align*}
&  \underbrace{\left(  -4\right)  ^{n}}_{=\left(  -1\right)  ^{n}\cdot4^{n}%
}\cdot\left(  -1\right)  ^{n/2}\underbrace{\dbinom{-1/2}{n/2}}_{=\dfrac
{\left(  -1\right)  ^{n/2}}{4^{n/2}}\dbinom{n}{n/2}}\\
&  =\left(  -1\right)  ^{n}\cdot4^{n}\cdot\left(  -1\right)  ^{n/2}%
\dfrac{\left(  -1\right)  ^{n/2}}{4^{n/2}}\dbinom{n}{n/2}\\
&  =\underbrace{\left(  -1\right)  ^{n}}_{\substack{=1\\\text{(since }n\text{
is even)}}}\cdot\underbrace{\left(  -1\right)  ^{n/2}\left(  -1\right)
^{n/2}}_{\substack{=\left(  -1\right)  ^{n/2+n/2}=1\\\text{(since
}n/2+n/2=n\text{ is even)}}}\cdot\underbrace{\dfrac{4^{n}}{4^{n/2}}%
}_{=4^{n-n/2}=2^{n}}\dbinom{n}{n/2}\\
&  =2^{n}\dbinom{n}{n/2}.
\end{align*}

Now, forget that we assumed that $n$ is even. We thus have proven that if $n$
is even, then $\left(  -4\right)  ^{n}\cdot\left(  -1\right)  ^{n/2}%
\dbinom{-1/2}{n/2}=2^{n}\dbinom{n}{n/2}$. Hence,%
\[%
\begin{cases}
\left(  -4\right)  ^{n}\cdot\left(  -1\right)  ^{n/2}\dbinom{-1/2}{n/2}, &
\text{if }n\text{ is even};\\
0, & \text{if }n\text{ is odd}%
\end{cases}
=%
\begin{cases}
2^{n}\dbinom{n}{n/2}, & \text{if }n\text{ is even};\\
0, & \text{if }n\text{ is odd}%
\end{cases}
.
\]
Thus, (\ref{sol.vander-1/2.b.step2}) becomes%
\begin{align*}
&  \sum_{k=0}^{n}\left(  -1\right)  ^{k}\dbinom{2k}{k}\dbinom{2\left(
n-k\right)  }{n-k}\\
&  =%
\begin{cases}
\left(  -4\right)  ^{n}\cdot\left(  -1\right)  ^{n/2}\dbinom{-1/2}{n/2}, &
\text{if }n\text{ is even};\\
0, & \text{if }n\text{ is odd}%
\end{cases}
=%
\begin{cases}
2^{n}\dbinom{n}{n/2}, & \text{if }n\text{ is even};\\
0, & \text{if }n\text{ is odd}%
\end{cases}
.
\end{align*}
This solves Exercise \ref{exe.vander-1/2} \textbf{(b)}.
\end{proof}

\subsection{Solution to Exercise \ref{exe.central-binomial-even}}

\begin{vershort}
\begin{proof}
[Solution to Exercise \ref{exe.central-binomial-even}.]Proposition
\ref{prop.binom.int} (applied to $2m$ and $m$ instead of $m$ and $n$) yields
$\dbinom{2m}{m}\in\mathbb{Z}$. Similarly, we can find $\dbinom{2m-1}{m-1}%
\in\mathbb{Z}$. Thus, it makes sense to speak of $\dbinom{2m}{m}$ or
$\dbinom{2m-1}{m-1}$ being even.

Proposition \ref{prop.binom.X-1} (applied to $2m$ and $m$ instead of $m$ and
$n$) yields $\dbinom{2m}{m}=\underbrace{\dfrac{2m}{m}}_{=2}\dbinom{2m-1}%
{m-1}=2\dbinom{2m-1}{m-1}$.

\textbf{(a)} We have $\dbinom{2m-1}{m-1}\in\mathbb{Z}$. Thus, the integer
$2\dbinom{2m-1}{m-1}$ is even. In view of $\dbinom{2m}{m}=2\dbinom{2m-1}{m-1}%
$, this rewrites as follows: The integer $\dbinom{2m}{m}$ is even. This solves
Exercise \ref{exe.central-binomial-even} \textbf{(a)}.

\textbf{(b)} Assume that $m$ is odd and satisfies $m>1$. We have
$m\equiv1\operatorname{mod}2$ (since $m$ is odd). Also, $m-1$ is a positive
integer (since $m>1$). Hence, Exercise \ref{exe.central-binomial-even}
\textbf{(a)} (applied to $m-1$ instead of $m$) shows that $\dbinom{2\left(
m-1\right)  }{m-1}$ is even. In other words, $\dbinom{2\left(  m-1\right)
}{m-1}\equiv0\operatorname{mod}2$.

But $2m-1=\underbrace{m}_{>1}+m-1>1+m-1=m$, so that $2m-1\geq m\geq0$ and thus
$2m-1\in\mathbb{N}$. Hence, Proposition \ref{prop.binom.symm} (applied to
$2m-1$ and $m$ instead of $m$ and $n$) yields $\dbinom{2m-1}{m}=\dbinom
{2m-1}{\left(  2m-1\right)  -m}=\dbinom{2m-1}{m-1}$ (since $\left(
2m-1\right)  -m=m-1$).

But Proposition \ref{prop.binom.X-1} (applied to $2m-1$ and $m$ instead of $m$
and $n$) yields%
\[
\dbinom{2m-1}{m}=\dfrac{2m-1}{m}\dbinom{\left(  2m-1\right)  -1}{m-1}%
=\dfrac{2m-1}{m}\dbinom{2\left(  m-1\right)  }{m-1}%
\]
(since $\left(  2m-1\right)  -1=2\left(  m-1\right)  $). Multiplying both
sides of this equality by $m$, we find%
\[
m\dbinom{2m-1}{m}=\left(  2m-1\right)  \underbrace{\dbinom{2\left(
m-1\right)  }{m-1}}_{\equiv0\operatorname{mod}2}\equiv0\operatorname{mod}2.
\]
Hence,%
\[
0\equiv\underbrace{m}_{\equiv1\operatorname{mod}2}\dbinom{2m-1}{m}%
\equiv\dbinom{2m-1}{m}=\dbinom{2m-1}{m-1}\operatorname{mod}2.
\]
In other words, $\dbinom{2m-1}{m-1}\equiv0\operatorname{mod}2$. In other
words, the integer $\dbinom{2m-1}{m-1}$ is even. This solves Exercise
\ref{exe.central-binomial-even} \textbf{(b)}.

\textbf{(c)} Assume that $m$ is odd and satisfies $m>1$. Exercise
\ref{exe.central-binomial-even} \textbf{(b)} yields that the integer
$\dbinom{2m-1}{m-1}$ is even. In other words, there exists an integer $z$ such
that $\dbinom{2m-1}{m-1}=2z$. Consider this $z$. Now,%
\[
\dbinom{2m}{m}=2\underbrace{\dbinom{2m-1}{m-1}}_{=2z}=2\cdot2z=4z\equiv
0\operatorname{mod}4
\]
(since $z$ is an integer). This solves Exercise
\ref{exe.central-binomial-even} \textbf{(c)}.
\end{proof}
\end{vershort}

\begin{verlong}
\begin{proof}
[Solution to Exercise \ref{exe.central-binomial-even}.]We have $m\in
\mathbb{N}$ (since $m$ is a positive integer), so that $2m\in\mathbb{N}%
\subseteq\mathbb{Z}$. Hence, Proposition \ref{prop.binom.int} (applied to $2m$
and $m$ instead of $m$ and $n$) yields $\dbinom{2m}{m}\in\mathbb{Z}$. Thus, it
makes sense to speak of $\dbinom{2m}{m}$ being even.

Also, $m-1\in\mathbb{N}$ (since $m$ is a positive integer) and $2m-1\in
\mathbb{Z}$ (since $m\in\mathbb{N}$). Hence, Proposition \ref{prop.binom.int}
(applied to $2m-1$ and $m-1$ instead of $m$ and $n$) yields $\dbinom
{2m-1}{m-1}\in\mathbb{Z}$. Thus, it makes sense to speak of $\dbinom
{2m-1}{m-1}$ being even.

We have $m\in\left\{  1,2,3,\ldots\right\}  $ (since $m$ is a positive
integer). Thus, Proposition \ref{prop.binom.X-1} (applied to $2m$ and $m$
instead of $m$ and $n$) yields $\dbinom{2m}{m}=\underbrace{\dfrac{2m}{m}}%
_{=2}\dbinom{2m-1}{m-1}=2\dbinom{2m-1}{m-1}$.

\textbf{(a)} We have $\dbinom{2m-1}{m-1}\in\mathbb{Z}$. Thus, the integer
$2\dbinom{2m-1}{m-1}$ is even. In view of $\dbinom{2m}{m}=2\dbinom{2m-1}{m-1}%
$, this rewrites as follows: The integer $\dbinom{2m}{m}$ is even. This solves
Exercise \ref{exe.central-binomial-even} \textbf{(a)}.

\textbf{(b)} Assume that $m$ is odd and satisfies $m>1$. We have
$m\equiv1\operatorname{mod}2$ (since $m$ is odd). Also, $m-1>0$ (since $m>1$);
thus, $m-1$ is a positive integer (since $m-1$ is an integer). Hence, Exercise
\ref{exe.central-binomial-even} \textbf{(a)} (applied to $m-1$ instead of $m$)
shows that $\dbinom{2\left(  m-1\right)  }{m-1}$ is even. In other words,
$\dbinom{2\left(  m-1\right)  }{m-1}\equiv0\operatorname{mod}2$.

But $2m-1=\underbrace{m}_{>1}+m-1>1+m-1=m$, so that $2m-1\geq m\geq0$ (since
$m\in\mathbb{N}$) and thus $2m-1\in\mathbb{N}$. Hence, Proposition
\ref{prop.binom.symm} (applied to $2m-1$ and $m$ instead of $m$ and $n$)
yields $\dbinom{2m-1}{m}=\dbinom{2m-1}{\left(  2m-1\right)  -m}=\dbinom
{2m-1}{m-1}$ (since $\left(  2m-1\right)  -m=m-1$).

But $m\in\left\{  1,2,3,\ldots\right\}  $ (since $m$ is a positive integer).
Hence, Proposition \ref{prop.binom.X-1} (applied to $2m-1$ and $m$ instead of
$m$ and $n$) yields%
\[
\dbinom{2m-1}{m}=\dfrac{2m-1}{m}\dbinom{\left(  2m-1\right)  -1}{m-1}.
\]
Multiplying both sides of this equality by $m$, we find%
\begin{align*}
m\dbinom{2m-1}{m}  &  =\left(  2m-1\right)  \dbinom{\left(  2m-1\right)
-1}{m-1}=\left(  2m-1\right)  \underbrace{\dbinom{2\left(  m-1\right)  }{m-1}%
}_{\equiv0\operatorname{mod}2}\\
&  \ \ \ \ \ \ \ \ \ \ \left(  \text{since }\left(  2m-1\right)  -1=2\left(
m-1\right)  \right) \\
&  \equiv0\operatorname{mod}2.
\end{align*}
Hence,%
\[
0\equiv\underbrace{m}_{\equiv1\operatorname{mod}2}\dbinom{2m-1}{m}%
\equiv\dbinom{2m-1}{m}=\dbinom{2m-1}{m-1}\operatorname{mod}2.
\]
In other words, $\dbinom{2m-1}{m-1}\equiv0\operatorname{mod}2$. In other
words, the integer $\dbinom{2m-1}{m-1}$ is even. This solves Exercise
\ref{exe.central-binomial-even} \textbf{(b)}.

\textbf{(c)} Assume that $m$ is odd and satisfies $m>1$. Exercise
\ref{exe.central-binomial-even} \textbf{(b)} yields that the integer
$\dbinom{2m-1}{m-1}$ is even. In other words, there exists an integer $z$ such
that $\dbinom{2m-1}{m-1}=2z$. Consider this $z$. Now,%
\[
\dbinom{2m}{m}=2\underbrace{\dbinom{2m-1}{m-1}}_{=2z}=2\cdot2z=4z\equiv
0\operatorname{mod}4
\]
(since $z$ is an integer). This solves Exercise
\ref{exe.central-binomial-even} \textbf{(c)}.
\end{proof}
\end{verlong}

\subsection{Solution to Exercise \ref{exe.supercat}}

Before we solve Exercise \ref{exe.supercat}, let us prove some straightforward identities:

\begin{lemma}
\label{lem.binom.symmetry-m+n}Let $m\in\mathbb{N}$ and $n\in\mathbb{N}$. Then,
$\dbinom{m+n}{m}=\dbinom{m+n}{n}$.
\end{lemma}

\begin{proof}
[Proof of Lemma \ref{lem.binom.symmetry-m+n}.]We have $m+n\in\mathbb{N}$
(since $m\in\mathbb{N}$ and $n\in\mathbb{N}$). Also, $n\in\mathbb{N}$, so that
$n\geq0$. Now, $m+\underbrace{n}_{\geq0}\geq m$. Hence, Proposition
\ref{prop.binom.symm} (applied to $m+n$ and $m$ instead of $m$ and $n$) yields
$\dbinom{m+n}{m}=\dbinom{m+n}{\left(  m+n\right)  -m}=\dbinom{m+n}{n}$ (since
$\left(  m+n\right)  -m=n$). This proves Lemma \ref{lem.binom.symmetry-m+n}.
\end{proof}

\begin{lemma}
\label{lem.sol.supercat.2}Let $m\in\mathbb{N}$ and $n\in\mathbb{N}$. Let
$p=\min\left\{  m,n\right\}  $. Let $k\in\left\{  -p,-p+1,\ldots,p\right\}  $.
Then,%
\[
\dbinom{2m}{m+k}\dbinom{2n}{n-k}=\dfrac{\left(  2m\right)  !\left(  2n\right)
!}{\left(  m+n\right)  !^{2}}\cdot\dbinom{m+n}{m+k}\dbinom{m+n}{n+k}.
\]

\end{lemma}

\begin{proof}
[Proof of Lemma \ref{lem.sol.supercat.2}.]Clearly, $p=\min\left\{
m,n\right\}  \in\mathbb{N}$ (since $m\in\mathbb{N}$ and $n\in\mathbb{N}$). We
have $k\in\left\{  -p,-p+1,\ldots,p\right\}  $. Thus, $k$ is an integer
satisfying $-p\leq k\leq p$.

We have $p=\min\left\{  m,n\right\}  \leq m$. Hence, $k\leq p\leq m$. Thus,
$m-k\geq0$. Also, $-p\leq k$, so that $k\geq-\underbrace{p}_{\leq m}\geq-m$.
Hence, $k+m\geq0$.

We have $p=\min\left\{  m,n\right\}  \leq n$. Hence, $k\leq p\leq n$. Thus,
$n-k\geq0$. Also, $-p\leq k$, so that $k\geq-\underbrace{p}_{\leq n}\geq-n$.
Hence, $k+n\geq0$.

\begin{vershort}
From $n-k\geq0$, we obtain $n-k\in\mathbb{N}$ (since $n-k$ is an integer).
Also, $m+n\in\mathbb{N}$ (since $m\in\mathbb{N}$ and $n\in\mathbb{N}$).
Furthermore, $n+k=k+n\geq0$. Hence, $n+k\in\mathbb{N}$ (since $n+k$ is an
integer). Similarly, $m+k\in\mathbb{N}$.
\end{vershort}

\begin{verlong}
From $n-k\geq0$, we obtain $n-k\in\mathbb{N}$ (since $n-k$ is an integer
(because $n$ and $k$ are integers)). Also, $m+n\in\mathbb{N}$ (since
$m\in\mathbb{N}$ and $n\in\mathbb{N}$). Furthermore, $n+k=k+n\geq0$. Hence,
$n+k\in\mathbb{N}$ (since $n+k$ is an integer (since $n$ and $k$ are
integers)). Also, $m+k=k+m\geq0$. Hence, $m+k\in\mathbb{N}$ (since $m+k$ is an
integer (since $m$ and $k$ are integers)).
\end{verlong}

\begin{verlong}
The number $\left(  m+n\right)  !=1\cdot2\cdot\cdots\cdot\left(  m+n\right)  $
is a positive integer. Hence, its square $\left(  m+n\right)  !^{2}$ is a
positive integer as well, and thus is nonzero. Hence, the fraction
$\dfrac{\left(  2m\right)  !\left(  2n\right)  !}{\left(  m+n\right)  !^{2}}$
is well-defined.
\end{verlong}

Now, $m+n\geq m+k$ (since $m+\underbrace{k}_{\leq n}\leq m+n$). Hence,
Proposition \ref{prop.binom.formula} (applied to $m+n$ and $m+k$ instead of
$m$ and $n$) yields%
\begin{equation}
\dbinom{m+n}{m+k}=\dfrac{\left(  m+n\right)  !}{\left(  m+k\right)  !\left(
\left(  m+n\right)  -\left(  m+k\right)  \right)  !}=\dfrac{\left(
m+n\right)  !}{\left(  m+k\right)  !\left(  n-k\right)  !}
\label{pf.lem.sol.supercat.2.R1}%
\end{equation}
(since $\left(  m+n\right)  -\left(  m+k\right)  =n-k$).

Also, $m+n\geq n+k$ (since $n+\underbrace{k}_{\leq m}\leq n+m=m+n$). Hence,
Proposition \ref{prop.binom.formula} (applied to $m+n$ and $n+k$ instead of
$m$ and $n$) yields%
\begin{equation}
\dbinom{m+n}{n+k}=\dfrac{\left(  m+n\right)  !}{\left(  n+k\right)  !\left(
\left(  m+n\right)  -\left(  n+k\right)  \right)  !}=\dfrac{\left(
m+n\right)  !}{\left(  n+k\right)  !\left(  m-k\right)  !}
\label{pf.lem.sol.supercat.2.R2}%
\end{equation}
(since $\left(  m+n\right)  -\left(  n+k\right)  =m-k$).

Now,%
\begin{align}
&  \dfrac{\left(  2m\right)  !\left(  2n\right)  !}{\left(  m+n\right)  !^{2}%
}\cdot\underbrace{\dbinom{m+n}{m+k}}_{\substack{=\dfrac{\left(  m+n\right)
!}{\left(  m+k\right)  !\left(  n-k\right)  !}\\\text{(by
(\ref{pf.lem.sol.supercat.2.R1}))}}}\ \ \underbrace{\dbinom{m+n}{n+k}%
}_{\substack{=\dfrac{\left(  m+n\right)  !}{\left(  n+k\right)  !\left(
m-k\right)  !}\\\text{(by (\ref{pf.lem.sol.supercat.2.R2}))}}}\nonumber\\
&  =\dfrac{\left(  2m\right)  !\left(  2n\right)  !}{\left(  m+n\right)
!^{2}}\cdot\dfrac{\left(  m+n\right)  !}{\left(  m+k\right)  !\left(
n-k\right)  !}\cdot\dfrac{\left(  m+n\right)  !}{\left(  n+k\right)  !\left(
m-k\right)  !}\nonumber\\
&  =\dfrac{\left(  2m\right)  !\left(  2n\right)  !}{\left(  m+k\right)
!\left(  n-k\right)  !\left(  n+k\right)  !\left(  m-k\right)  !}.
\label{pf.lem.sol.supercat.2.R}%
\end{align}

Also, $2m\in\mathbb{N}$ (since $m\in\mathbb{N}$) and $2m\geq m+k$ (since
$2m-\left(  m+k\right)  =m-k\geq0$). Hence, Proposition
\ref{prop.binom.formula} (applied to $2m$ and $m+k$ instead of $m$ and $n$)
yields%
\begin{equation}
\dbinom{2m}{m+k}=\dfrac{\left(  2m\right)  !}{\left(  m+k\right)  !\left(
2m-\left(  m+k\right)  \right)  !}=\dfrac{\left(  2m\right)  !}{\left(
m+k\right)  !\left(  m-k\right)  !} \label{pf.lem.sol.supercat.2.L1}%
\end{equation}
(since $2m-\left(  m+k\right)  =m-k$).

Moreover, $2n\in\mathbb{N}$ (since $n\in\mathbb{N}$) and $2n\geq n-k$ (since
$2n-\left(  n-k\right)  =k+n\geq0$). Hence, Proposition
\ref{prop.binom.formula} (applied to $2n$ and $n-k$ instead of $m$ and $n$)
yields%
\begin{equation}
\dbinom{2n}{n-k}=\dfrac{\left(  2n\right)  !}{\left(  n-k\right)  !\left(
2n-\left(  n-k\right)  \right)  !}=\dfrac{\left(  2n\right)  !}{\left(
n-k\right)  !\left(  n+k\right)  !} \label{pf.lem.sol.supercat.2.L2}%
\end{equation}
(since $2n-\left(  n-k\right)  =n+k$).

Multiplying the equalities (\ref{pf.lem.sol.supercat.2.L1}) and
(\ref{pf.lem.sol.supercat.2.L2}), we obtain%
\begin{align*}
\dbinom{2m}{m+k}\dbinom{2n}{n-k}  &  =\dfrac{\left(  2m\right)  !}{\left(
m+k\right)  !\left(  m-k\right)  !}\cdot\dfrac{\left(  2n\right)  !}{\left(
n-k\right)  !\left(  n+k\right)  !}\\
&  =\dfrac{\left(  2m\right)  !\left(  2n\right)  !}{\left(  m+k\right)
!\left(  n-k\right)  !\left(  n+k\right)  !\left(  m-k\right)  !}\\
&  =\dfrac{\left(  2m\right)  !\left(  2n\right)  !}{\left(  m+n\right)
!^{2}}\cdot\dbinom{m+n}{m+k}\dbinom{m+n}{n+k}%
\end{align*}
(by (\ref{pf.lem.sol.supercat.2.R})). This proves Lemma
\ref{lem.sol.supercat.2}.
\end{proof}

\begin{proof}
[Solution to Exercise \ref{exe.supercat}.]Let us begin with parts
\textbf{(e)}, \textbf{(f)}, \textbf{(g)}, \textbf{(h)} and \textbf{(i)} of
Exercise \ref{exe.supercat}.

\begin{vershort}
\textbf{(e)} Let $m\in\mathbb{N}$. Thus, $2m\geq m\geq0$. Hence, Proposition
\ref{prop.binom.formula} (applied to $2m$ and $m$ instead of $m$ and $n$)
yields $\dbinom{2m}{m}=\dfrac{\left(  2m\right)  !}{m!\left(  2m-m\right)
!}=\dfrac{\left(  2m\right)  !}{m!m!}$. But the definition of $T\left(
m,0\right)  $ yields%
\begin{align*}
T\left(  m,0\right)   &  =\dfrac{\left(  2m\right)  !\left(  2\cdot0\right)
!}{m!0!\left(  m+0\right)  !}=\dfrac{\left(  2m\right)  !\cdot1}{m!\cdot
1\cdot\left(  m+0\right)  !}\ \ \ \ \ \ \ \ \ \ \left(  \text{since }\left(
2\cdot0\right)  !=0!=1\text{ and }0!=1\right) \\
&  =\dfrac{\left(  2m\right)  !}{m!\left(  m+0\right)  !}=\dfrac{\left(
2m\right)  !}{m!m!}=\dbinom{2m}{m}.
\end{align*}
This solves Exercise \ref{exe.supercat} \textbf{(e)}.
\end{vershort}

\begin{verlong}
\textbf{(e)} Let $m\in\mathbb{N}$. Thus, $m\geq0$, so that $2m\geq m$ (since
$2m-m=m\geq0$) and therefore $2m\geq m\geq0$. Thus, $2m\in\mathbb{N}$. Hence,
Proposition \ref{prop.binom.formula} (applied to $2m$ and $m$ instead of $m$
and $n$) yields $\dbinom{2m}{m}=\dfrac{\left(  2m\right)  !}{m!\left(
2m-m\right)  !}=\dfrac{\left(  2m\right)  !}{m!m!}$ (since $2m-m=m$).

But the definition of $T\left(  m,0\right)  $ yields%
\begin{align*}
T\left(  m,0\right)   &  =\dfrac{\left(  2m\right)  !\left(  2\cdot0\right)
!}{m!0!\left(  m+0\right)  !}=\dfrac{\left(  2m\right)  !\cdot1}{m!\cdot
1\cdot\left(  m+0\right)  !}\ \ \ \ \ \ \ \ \ \ \left(  \text{since }\left(
\underbrace{2\cdot0}_{=0}\right)  !=0!=1\text{ and }0!=1\right) \\
&  =\dfrac{\left(  2m\right)  !}{m!\left(  m+0\right)  !}=\dfrac{\left(
2m\right)  !}{m!m!}\ \ \ \ \ \ \ \ \ \ \left(  \text{since }m+0=m\right) \\
&  =\dbinom{2m}{m}.
\end{align*}
This solves Exercise \ref{exe.supercat} \textbf{(e)}.
\end{verlong}

\textbf{(f)} Let $m\in\mathbb{N}$ and $n\in\mathbb{N}$.

\begin{vershort}
From $m\in\mathbb{N}$, we obtain $2m\geq m$. Hence, Proposition
\ref{prop.binom.formula} (applied to $2m$ and $m$ instead of $m$ and $n$)
yields $\dbinom{2m}{m}=\dfrac{\left(  2m\right)  !}{m!\left(  2m-m\right)
!}=\dfrac{\left(  2m\right)  !}{m!m!}$. Similarly, $\dbinom{2n}{n}%
=\dfrac{\left(  2n\right)  !}{n!n!}$. Also, $m+\underbrace{n}_{\geq0}\geq m$.
Hence, Proposition \ref{prop.binom.formula} (applied to $m+n$ and $m$ instead
of $m$ and $n$) yields $\dbinom{m+n}{m}=\dfrac{\left(  m+n\right)
!}{m!\left(  \left(  m+n\right)  -m\right)  !}=\dfrac{\left(  m+n\right)
!}{m!n!}$. Thus, the rational number $\dbinom{m+n}{m}$ is nonzero (since
$\dfrac{\left(  m+n\right)  !}{m!n!}$ is clearly nonzero). Hence, the fraction
$\dfrac{\dbinom{2m}{m}\dbinom{2n}{n}}{\dbinom{m+n}{m}}$ is well-defined.
\end{vershort}

\begin{verlong}
From $m\in\mathbb{N}$, we obtain $m\geq0$, so that $2m\geq m$ (since
$2m-m=m\geq0$) and therefore $2m\geq m\geq0$. Thus, $2m\in\mathbb{N}$. Hence,
Proposition \ref{prop.binom.formula} (applied to $2m$ and $m$ instead of $m$
and $n$) yields $\dbinom{2m}{m}=\dfrac{\left(  2m\right)  !}{m!\left(
2m-m\right)  !}=\dfrac{\left(  2m\right)  !}{m!m!}$ (since $2m-m=m$). The same
argument (applied to $n$ instead of $m$) yields $\dbinom{2n}{n}=\dfrac{\left(
2n\right)  !}{n!n!}$.

Clearly, $\left(  m+n\right)  !=1\cdot2\cdot\cdots\cdot\left(  m+n\right)  $
is a positive integer, and thus a nonzero rational number.

Also, $m+n\in\mathbb{N}$ (since $m\in\mathbb{N}$ and $n\in\mathbb{N}$) and
$m+\underbrace{n}_{\geq0}\geq m$. Hence, Proposition \ref{prop.binom.formula}
(applied to $m+n$ and $m$ instead of $m$ and $n$) yields $\dbinom{m+n}%
{m}=\dfrac{\left(  m+n\right)  !}{m!\left(  \left(  m+n\right)  -m\right)
!}=\dfrac{\left(  m+n\right)  !}{m!n!}$ (since $\left(  m+n\right)  -m=n$).
Thus, the rational number $\dbinom{m+n}{m}$ is nonzero (since $\dfrac{\left(
m+n\right)  !}{m!n!}$ is nonzero (because $\left(  m+n\right)  !$ is
nonzero)). Hence, the fraction $\dfrac{\dbinom{2m}{m}\dbinom{2n}{n}}%
{\dbinom{m+n}{m}}$ is well-defined.
\end{verlong}

Now,%
\begin{align}
\dfrac{\dbinom{2m}{m}\dbinom{2n}{n}}{\dbinom{m+n}{m}}  &  =\dfrac
{\dfrac{\left(  2m\right)  !}{m!m!}\cdot\dfrac{\left(  2n\right)  !}{n!n!}%
}{\left(  \dfrac{\left(  m+n\right)  !}{m!n!}\right)  }%
\ \ \ \ \ \ \ \ \ \ \left(
\begin{array}
[c]{c}%
\text{since }\dbinom{2m}{m}=\dfrac{\left(  2m\right)  !}{m!m!}\text{ and
}\dbinom{2n}{n}=\dfrac{\left(  2n\right)  !}{n!n!}\\
\text{and }\dbinom{m+n}{m}=\dfrac{\left(  m+n\right)  !}{m!n!}%
\end{array}
\right) \nonumber\\
&  =\dfrac{\left(  2m\right)  !\left(  2n\right)  !}{m!n!\left(  m+n\right)
!}. \label{sol.supercat.f.1}%
\end{align}
But the definition of $T\left(  m,n\right)  $ yields
\[
T\left(  m,n\right)  =\dfrac{\left(  2m\right)  !\left(  2n\right)
!}{m!n!\left(  m+n\right)  !}.
\]
Comparing this with (\ref{sol.supercat.f.1}), we obtain $T\left(  m,n\right)
=\dfrac{\dbinom{2m}{m}\dbinom{2n}{n}}{\dbinom{m+n}{m}}$. This solves Exercise
\ref{exe.supercat} \textbf{(f)}.

\textbf{(g)} Let $m\in\mathbb{N}$ and $n\in\mathbb{N}$. The definition of
$T\left(  n,m\right)  $ yields%
\begin{equation}
T\left(  n,m\right)  =\dfrac{\left(  2n\right)  !\left(  2m\right)
!}{n!m!\left(  n+m\right)  !}=\dfrac{\left(  2m\right)  !\left(  2n\right)
!}{m!n!\left(  n+m\right)  !}=\dfrac{\left(  2m\right)  !\left(  2n\right)
!}{m!n!\left(  m+n\right)  !} \label{sol.supercat.g.1}%
\end{equation}
(since $n+m=m+n$). But the definition of $T\left(  m,n\right)  $ yields
\[
T\left(  m,n\right)  =\dfrac{\left(  2m\right)  !\left(  2n\right)
!}{m!n!\left(  m+n\right)  !}.
\]
Comparing this with (\ref{sol.supercat.g.1}), we obtain $T\left(  m,n\right)
=T\left(  n,m\right)  $. This solves Exercise \ref{exe.supercat} \textbf{(g)}.

\textbf{(h)} Let us forget that we fixed $m$, $n$ and $p$. We shall first
prove some auxiliary observations:

\begin{statement}
\textit{Observation 1:} Let $m\in\mathbb{N}$ and $n\in\mathbb{N}$ be such that
$m\geq n$. Then,
\[
\sum_{k=-n}^{n}\left(  -1\right)  ^{k}\dbinom{m+n}{m+k}\dbinom{m+n}%
{n+k}=\dbinom{m+n}{m}.
\]

\end{statement}

\begin{vershort}
[\textit{Proof of Observation 1:} Every $k\in\left\{  -n,-n+1,\ldots
,n\right\}  $ satisfies%
\begin{equation}
\dbinom{m+n}{n-k}=\dbinom{m+n}{m+k} \label{sol.supercat.h.o1.pf.short.1}%
\end{equation}
\footnote{\textit{Proof of (\ref{sol.supercat.h.o1.pf.short.1}):} Let
$k\in\left\{  -n,-n+1,\ldots,n\right\}  $. Thus, $k$ is an integer satisfying
$-n\leq k\leq n$. From $k\leq n$, we obtain $n-k\geq0$. Hence, $n-k\in
\mathbb{N}$ (since $n-k$ is an integer). Also, $\left(  m+n\right)  -\left(
n-k\right)  =\underbrace{m}_{\geq n}+\underbrace{k}_{\substack{\geq
-n\\\text{(since }-n\leq k\text{)}}}\geq n+\left(  -n\right)  =0$; in other
words, $m+n\geq n-k$. Hence, Proposition \ref{prop.binom.symm} (applied to
$m+n$ and $n-k$ instead of $m$ and $n$) yields $\dbinom{m+n}{n-k}=\dbinom
{m+n}{\left(  m+n\right)  -\left(  n-k\right)  }=\dbinom{m+n}{m+k}$. This
proves (\ref{sol.supercat.h.o1.pf.short.1}).}.

Exercise \ref{exe.AoPS262752} (applied to $2n$ instead of $n$) yields%
\begin{align*}
\sum_{k=0}^{2n}\left(  -1\right)  ^{k}\dbinom{X}{k}\dbinom{X}{2n-k}  &  =%
\begin{cases}
\left(  -1\right)  ^{2n/2}\dbinom{X}{2n/2}, & \text{if }2n\text{ is even};\\
0, & \text{if }2n\text{ is odd}%
\end{cases}
\\
&  =\left(  -1\right)  ^{2n/2}\dbinom{X}{2n/2}\ \ \ \ \ \ \ \ \ \ \left(
\text{since }2n\text{ is even}\right) \\
&  =\left(  -1\right)  ^{n}\dbinom{X}{n}\ \ \ \ \ \ \ \ \ \ \left(
\text{since }2n/2=n\right)
\end{align*}
(an identity between polynomials in $\mathbb{Q}\left[  X\right]  $).
Substituting $m+n$ for $X$ in this equality, we obtain%
\[
\sum_{k=0}^{2n}\left(  -1\right)  ^{k}\dbinom{m+n}{k}\dbinom{m+n}%
{2n-k}=\left(  -1\right)  ^{n}\dbinom{m+n}{n}.
\]
Hence,%
\begin{align*}
\left(  -1\right)  ^{n}\dbinom{m+n}{n}  &  =\sum_{k=0}^{2n}\left(  -1\right)
^{k}\dbinom{m+n}{k}\dbinom{m+n}{2n-k}\\
&  =\sum_{k=-n}^{n}\underbrace{\left(  -1\right)  ^{k+n}}_{=\left(  -1\right)
^{k}\left(  -1\right)  ^{n}}\underbrace{\dbinom{m+n}{k+n}}_{=\dbinom{m+n}%
{n+k}}\underbrace{\dbinom{m+n}{2n-\left(  k+n\right)  }}_{\substack{=\dbinom
{m+n}{n-k}\\\text{(since }2n-\left(  k+n\right)  =n-k\text{)}}}\\
&  \ \ \ \ \ \ \ \ \ \ \left(  \text{here, we have substituted }k+n\text{ for
}k\text{ in the sum}\right) \\
&  =\left(  -1\right)  ^{n}\sum_{k=-n}^{n}\left(  -1\right)  ^{k}\dbinom
{m+n}{n+k}\dbinom{m+n}{n-k}.
\end{align*}
Dividing both sides of this equality by $\left(  -1\right)  ^{n}$, we obtain%
\begin{align*}
\dbinom{m+n}{n}  &  =\sum_{k=-n}^{n}\left(  -1\right)  ^{k}\dbinom{m+n}%
{n+k}\underbrace{\dbinom{m+n}{n-k}}_{\substack{=\dbinom{m+n}{m+k}\\\text{(by
(\ref{sol.supercat.h.o1.pf.short.1}))}}}=\sum_{k=-n}^{n}\left(  -1\right)
^{k}\underbrace{\dbinom{m+n}{n+k}\dbinom{m+n}{m+k}}_{=\dbinom{m+n}{m+k}%
\dbinom{m+n}{n+k}}\\
&  =\sum_{k=-n}^{n}\left(  -1\right)  ^{k}\dbinom{m+n}{m+k}\dbinom{m+n}{n+k}.
\end{align*}
In other words,%
\[
\sum_{k=-n}^{n}\left(  -1\right)  ^{k}\dbinom{m+n}{m+k}\dbinom{m+n}%
{n+k}=\dbinom{m+n}{n}=\dbinom{m+n}{m}%
\]
(by Lemma \ref{lem.binom.symmetry-m+n}). This proves Observation 1.]
\end{vershort}

\begin{verlong}
[\textit{Proof of Observation 1:} Every $k\in\left\{  -n,-n+1,\ldots
,n\right\}  $ satisfies%
\begin{equation}
\dbinom{m+n}{n-k}=\dbinom{m+n}{m+k} \label{sol.supercat.h.o1.pf.1}%
\end{equation}
\footnote{\textit{Proof of (\ref{sol.supercat.h.o1.pf.1}):} Let $k\in\left\{
-n,-n+1,\ldots,n\right\}  $. Thus, $k$ is an integer satisfying $-n\leq k\leq
n$. From $k\leq n$, we obtain $n-k\geq0$. Hence, $n-k\in\mathbb{N}$ (since
$n-k$ is an integer (because $n$ and $k$ are integers)). Also, $m+n\in
\mathbb{N}$ (since $m\in\mathbb{N}$ and $n\in\mathbb{N}$). Finally, $\left(
m+n\right)  -\left(  n-k\right)  =\underbrace{m}_{\geq n}+\underbrace{k}%
_{\substack{\geq-n\\\text{(since }-n\leq k\text{)}}}\geq n+\left(  -n\right)
=0$; in other words, $m+n\geq n-k$. Hence, Proposition \ref{prop.binom.symm}
(applied to $m+n$ and $n-k$ instead of $m$ and $n$) yields $\dbinom{m+n}%
{n-k}=\dbinom{m+n}{\left(  m+n\right)  -\left(  n-k\right)  }=\dbinom
{m+n}{m+k}$ (since $\left(  m+n\right)  -\left(  n-k\right)  =m+k$). This
proves (\ref{sol.supercat.h.o1.pf.1}).}.

We have $2n\in\mathbb{N}$ (since $n\in\mathbb{N}$). Thus, Exercise
\ref{exe.AoPS262752} (applied to $2n$ instead of $n$) yields%
\begin{align*}
\sum_{k=0}^{2n}\left(  -1\right)  ^{k}\dbinom{X}{k}\dbinom{X}{2n-k}  &  =%
\begin{cases}
\left(  -1\right)  ^{2n/2}\dbinom{X}{2n/2}, & \text{if }2n\text{ is even};\\
0, & \text{if }2n\text{ is odd}%
\end{cases}
\\
&  =\left(  -1\right)  ^{2n/2}\dbinom{X}{2n/2}\ \ \ \ \ \ \ \ \ \ \left(
\text{since }2n\text{ is even}\right) \\
&  =\left(  -1\right)  ^{n}\dbinom{X}{n}\ \ \ \ \ \ \ \ \ \ \left(
\text{since }2n/2=n\right)
\end{align*}
(an identity between polynomials in $\mathbb{Q}\left[  X\right]  $).
Substituting $m+n$ for $X$ in this equality, we obtain%
\[
\sum_{k=0}^{2n}\left(  -1\right)  ^{k}\dbinom{m+n}{k}\dbinom{m+n}%
{2n-k}=\left(  -1\right)  ^{n}\dbinom{m+n}{n}.
\]
Hence,%
\begin{align*}
&  \left(  -1\right)  ^{n}\dbinom{m+n}{n}\\
&  =\sum_{k=0}^{2n}\left(  -1\right)  ^{k}\dbinom{m+n}{k}\dbinom{m+n}{2n-k}\\
&  =\underbrace{\sum_{k=-n}^{2n-n}}_{\substack{=\sum_{k=-n}^{n}\\\text{(since
}2n-n=n\text{)}}}\underbrace{\left(  -1\right)  ^{k+n}}_{=\left(  -1\right)
^{k}\left(  -1\right)  ^{n}}\underbrace{\dbinom{m+n}{k+n}}_{\substack{=\dbinom
{m+n}{n+k}\\\text{(since }k+n=n+k\text{)}}}\underbrace{\dbinom{m+n}{2n-\left(
k+n\right)  }}_{\substack{=\dbinom{m+n}{n-k}\\\text{(since }2n-\left(
k+n\right)  =n-k\text{)}}}\\
&  \ \ \ \ \ \ \ \ \ \ \left(  \text{here, we have substituted }k+n\text{ for
}k\text{ in the sum}\right) \\
&  =\sum_{k=-n}^{n}\left(  -1\right)  ^{k}\left(  -1\right)  ^{n}\dbinom
{m+n}{n+k}\dbinom{m+n}{n-k}\\
&  =\left(  -1\right)  ^{n}\sum_{k=-n}^{n}\left(  -1\right)  ^{k}\dbinom
{m+n}{n+k}\dbinom{m+n}{n-k}.
\end{align*}
We can divide both sides of this equality by $\left(  -1\right)  ^{n}$ (since
$\left(  -1\right)  ^{n}$ is a nonzero integer (because $\left(  -1\right)
^{n}\in\left\{  1,-1\right\}  $)). Thus, we obtain%
\begin{align*}
\dbinom{m+n}{n}  &  =\sum_{k=-n}^{n}\left(  -1\right)  ^{k}\dbinom{m+n}%
{n+k}\underbrace{\dbinom{m+n}{n-k}}_{\substack{=\dbinom{m+n}{m+k}\\\text{(by
(\ref{sol.supercat.h.o1.pf.1}))}}}\\
&  =\sum_{k=-n}^{n}\left(  -1\right)  ^{k}\underbrace{\dbinom{m+n}{n+k}%
\dbinom{m+n}{m+k}}_{=\dbinom{m+n}{m+k}\dbinom{m+n}{n+k}}\\
&  =\sum_{k=-n}^{n}\left(  -1\right)  ^{k}\dbinom{m+n}{m+k}\dbinom{m+n}{n+k}.
\end{align*}
In other words,%
\[
\sum_{k=-n}^{n}\left(  -1\right)  ^{k}\dbinom{m+n}{m+k}\dbinom{m+n}%
{n+k}=\dbinom{m+n}{n}=\dbinom{m+n}{m}%
\]
(since $\dbinom{m+n}{m}=\dbinom{m+n}{n}$ (by Lemma
\ref{lem.binom.symmetry-m+n})). This proves Observation 1.]
\end{verlong}

\begin{statement}
\textit{Observation 2:} Let $m\in\mathbb{N}$ and $n\in\mathbb{N}$ be such that
$m\leq n$. Then,
\[
\sum_{k=-m}^{m}\left(  -1\right)  ^{k}\dbinom{m+n}{m+k}\dbinom{m+n}%
{n+k}=\dbinom{m+n}{m}.
\]

\end{statement}

[\textit{Proof of Observation 2:} We have $m\leq n$, and thus $n\geq m$.
Hence, Observation 1 (applied to $n$ and $m$ instead of $m$ and $n$) yields
\[
\sum_{k=-m}^{m}\left(  -1\right)  ^{k}\dbinom{n+m}{n+k}\dbinom{n+m}%
{m+k}=\dbinom{n+m}{n}.
\]
Comparing this with%
\[
\sum_{k=-m}^{m}\left(  -1\right)  ^{k}\underbrace{\dbinom{n+m}{n+k}%
\dbinom{n+m}{m+k}}_{=\dbinom{n+m}{m+k}\dbinom{n+m}{n+k}}=\sum_{k=-m}%
^{m}\left(  -1\right)  ^{k}\dbinom{n+m}{m+k}\dbinom{n+m}{n+k},
\]
we obtain%
\[
\sum_{k=-m}^{m}\left(  -1\right)  ^{k}\dbinom{n+m}{m+k}\dbinom{n+m}%
{n+k}=\dbinom{n+m}{n}.
\]
In view of $n+m=m+n$, this rewrites as%
\[
\sum_{k=-m}^{m}\left(  -1\right)  ^{k}\dbinom{m+n}{m+k}\dbinom{m+n}%
{n+k}=\dbinom{m+n}{n}.
\]
Comparing this with
\[
\dbinom{m+n}{m}=\dbinom{m+n}{n}\ \ \ \ \ \ \ \ \ \ \left(  \text{by Lemma
\ref{lem.binom.symmetry-m+n}}\right)  ,
\]
we obtain%
\[
\sum_{k=-m}^{m}\left(  -1\right)  ^{k}\dbinom{m+n}{m+k}\dbinom{m+n}%
{n+k}=\dbinom{m+n}{m}.
\]
This proves Observation 2.]

Let us now come back to the solution of Exercise \ref{exe.supercat}
\textbf{(h)}. Let $m\in\mathbb{N}$ and $n\in\mathbb{N}$. Let $p=\min\left\{
m,n\right\}  $. We are in one of the following two cases:

\textit{Case 1:} We have $m\geq n$.

\textit{Case 2:} We have $m<n$.

Let us first consider Case 1. In this case, we have $m\geq n$. Now,
$p=\min\left\{  m,n\right\}  =n$ (since $m\geq n$). Hence,%
\[
\sum_{k=-p}^{p}\left(  -1\right)  ^{k}\dbinom{m+n}{m+k}\dbinom{m+n}{n+k}%
=\sum_{k=-n}^{n}\left(  -1\right)  ^{k}\dbinom{m+n}{m+k}\dbinom{m+n}%
{n+k}=\dbinom{m+n}{m}%
\]
(by Observation 1). Hence, Exercise \ref{exe.supercat} \textbf{(h)} is proven
in Case 1.

Let us first consider Case 2. In this case, we have $m<n$. Hence, $m\leq n$.
Now, $p=\min\left\{  m,n\right\}  =m$ (since $m\leq n$). Hence,%
\[
\sum_{k=-p}^{p}\left(  -1\right)  ^{k}\dbinom{m+n}{m+k}\dbinom{m+n}{n+k}%
=\sum_{k=-m}^{m}\left(  -1\right)  ^{k}\dbinom{m+n}{m+k}\dbinom{m+n}%
{n+k}=\dbinom{m+n}{m}%
\]
(by Observation 2). Hence, Exercise \ref{exe.supercat} \textbf{(h)} is proven
in Case 2.

\begin{vershort}
We have thus proven Exercise \ref{exe.supercat} \textbf{(h)} in each of the
two Cases 1 and 2. Thus, Exercise \ref{exe.supercat} \textbf{(h)} is solved in
all cases.
\end{vershort}

\begin{verlong}
We have thus proven Exercise \ref{exe.supercat} \textbf{(h)} in each of the
two Cases 1 and 2. Since these two Cases cover all possibilities, we thus
conclude that Exercise \ref{exe.supercat} \textbf{(h)} always holds. This
solves Exercise \ref{exe.supercat} \textbf{(h)}.
\end{verlong}

\begin{vershort}
\textbf{(i)} We have $m+n\geq m$ (since $n\geq0$). Hence, Proposition
\ref{prop.binom.formula} (applied to $m+n$ and $m$ instead of $m$ and $n$)
yields $\dbinom{m+n}{m}=\dfrac{\left(  m+n\right)  !}{m!\left(  \left(
m+n\right)  -m\right)  !}=\dfrac{\left(  m+n\right)  !}{m!n!}$.
\end{vershort}

\begin{verlong}
\textbf{(i)} We have $m+n\in\mathbb{N}$ (since $m\in\mathbb{N}$ and
$n\in\mathbb{N}$) and $m+\underbrace{n}_{\substack{\geq0\\\text{(since }%
n\in\mathbb{N}\text{)}}}\geq m$. Hence, Proposition \ref{prop.binom.formula}
(applied to $m+n$ and $m$ instead of $m$ and $n$) yields $\dbinom{m+n}%
{m}=\dfrac{\left(  m+n\right)  !}{m!\left(  \left(  m+n\right)  -m\right)
!}=\dfrac{\left(  m+n\right)  !}{m!n!}$ (since $\left(  m+n\right)  -m=n$).
\end{verlong}

But Exercise \ref{exe.supercat} \textbf{(h)} yields%
\begin{equation}
\sum_{k=-p}^{p}\left(  -1\right)  ^{k}\dbinom{m+n}{m+k}\dbinom{m+n}%
{n+k}=\dbinom{m+n}{m}=\dfrac{\left(  m+n\right)  !}{m!n!}.
\label{sol.supercat.i.o3.pf.2}%
\end{equation}

Now,
\begin{align*}
&  \sum_{k=-p}^{p}\left(  -1\right)  ^{k}\underbrace{\dbinom{2m}{m+k}%
\dbinom{2n}{n-k}}_{\substack{=\dfrac{\left(  2m\right)  !\left(  2n\right)
!}{\left(  m+n\right)  !^{2}}\cdot\dbinom{m+n}{m+k}\dbinom{m+n}{n+k}%
\\\text{(by Lemma \ref{lem.sol.supercat.2})}}}\\
&  =\sum_{k=-p}^{p}\left(  -1\right)  ^{k}\dfrac{\left(  2m\right)  !\left(
2n\right)  !}{\left(  m+n\right)  !^{2}}\cdot\dbinom{m+n}{m+k}\dbinom
{m+n}{n+k}\\
&  =\dfrac{\left(  2m\right)  !\left(  2n\right)  !}{\left(  m+n\right)
!^{2}}\cdot\underbrace{\sum_{k=-p}^{p}\left(  -1\right)  ^{k}\dbinom{m+n}%
{m+k}\dbinom{m+n}{n+k}}_{\substack{=\dfrac{\left(  m+n\right)  !}%
{m!n!}\\\text{(by (\ref{sol.supercat.i.o3.pf.2}))}}}\\
&  =\dfrac{\left(  2m\right)  !\left(  2n\right)  !}{\left(  m+n\right)
!^{2}}\cdot\dfrac{\left(  m+n\right)  !}{m!n!}=\dfrac{\left(  2m\right)
!\left(  2n\right)  !}{m!n!\left(  m+n\right)  !}.
\end{align*}
Comparing this with%
\[
T\left(  m,n\right)  =\dfrac{\left(  2m\right)  !\left(  2n\right)
!}{m!n!\left(  m+n\right)  !}\ \ \ \ \ \ \ \ \ \ \left(  \text{by the
definition of }T\left(  m,n\right)  \right)  ,
\]
we obtain $T\left(  m,n\right)  =\sum_{k=-p}^{p}\left(  -1\right)  ^{k}%
\dbinom{2m}{m+k}\dbinom{2n}{n-k}$. This solves Exercise \ref{exe.supercat}
\textbf{(i)}.

\begin{vershort}
\textbf{(a)} Let $m\in\mathbb{N}$ and $n\in\mathbb{N}$. Then, $2\left(
m+1\right)  =2m+2$, so that%
\begin{align}
\left(  2\left(  m+1\right)  \right)  !  &  =\left(  2m+2\right)
!=1\cdot2\cdot\cdots\cdot\left(  2m+2\right) \nonumber\\
&  =\underbrace{\left(  1\cdot2\cdot\cdots\cdot\left(  2m\right)  \right)
}_{=\left(  2m\right)  !}\cdot\left(  2m+1\right)  \cdot\underbrace{\left(
2m+2\right)  }_{=2\left(  m+1\right)  }\nonumber\\
&  =\left(  2m\right)  !\cdot\left(  2m+1\right)  \cdot2\left(  m+1\right)  .
\label{sol.supercat.a.short.1}%
\end{align}
Now, the definition of $T\left(  m+1,n\right)  $ yields%
\begin{align}
T\left(  m+1,n\right)   &  =\dfrac{\left(  2\left(  m+1\right)  \right)
!\left(  2n\right)  !}{\left(  m+1\right)  !n!\left(  m+1+n\right)  !}%
=\dfrac{\left(  2m\right)  !\cdot\left(  2m+1\right)  \cdot2\left(
m+1\right)  \cdot\left(  2n\right)  !}{m!\cdot\left(  m+1\right)  \cdot
n!\cdot\left(  m+n\right)  !\cdot\left(  m+n+1\right)  }\nonumber\\
&  \ \ \ \ \ \ \ \ \ \ \left(
\begin{array}
[c]{c}%
\text{by (\ref{sol.supercat.a.short.1}) and because of }\left(  m+1\right)
!=m!\cdot\left(  m+1\right) \\
\text{and }\left(  m+1+n\right)  !=\left(  m+n+1\right)  !=\left(  m+n\right)
!\cdot\left(  m+n+1\right)
\end{array}
\right) \nonumber\\
&  =\dfrac{2\cdot\left(  2m+1\right)  }{m+n+1}\cdot\underbrace{\dfrac{\left(
2m\right)  !\left(  2n\right)  !}{m!n!\left(  m+n\right)  !}}%
_{\substack{=T\left(  m,n\right)  \\\text{(by the definition of }T\left(
m,n\right)  \text{)}}}\nonumber\\
&  =\dfrac{2\cdot\left(  2m+1\right)  }{m+n+1}\cdot T\left(  m,n\right)  .
\label{sol.supercat.a.short.5}%
\end{align}
The same argument (but with $m$ and $n$ replaced by $n$ and $m$) shows that%
\[
T\left(  n+1,m\right)  =\dfrac{2\cdot\left(  2n+1\right)  }{n+m+1}\cdot
T\left(  n,m\right)  .
\]

But Exercise \ref{exe.supercat} \textbf{(g)} yields $T\left(  m,n\right)
=T\left(  n,m\right)  $. Also, Exercise \ref{exe.supercat} \textbf{(g)}
(applied to $n+1$ instead of $n$) yields%
\[
T\left(  m,n+1\right)  =T\left(  n+1,m\right)  =\underbrace{\dfrac
{2\cdot\left(  2n+1\right)  }{n+m+1}}_{=\dfrac{2\cdot\left(  2n+1\right)
}{m+n+1}}\cdot\underbrace{T\left(  n,m\right)  }_{=T\left(  m,n\right)
}=\dfrac{2\cdot\left(  2n+1\right)  }{m+n+1}\cdot T\left(  m,n\right)  .
\]
Adding this equality to (\ref{sol.supercat.a.short.5}), we obtain%
\begin{align*}
T\left(  m+1,n\right)  +T\left(  m,n+1\right)   &  =\dfrac{2\cdot\left(
2m+1\right)  }{m+n+1}\cdot T\left(  m,n\right)  +\dfrac{2\cdot\left(
2n+1\right)  }{m+n+1}\cdot T\left(  m,n\right) \\
&  =\underbrace{\left(  \dfrac{2\cdot\left(  2m+1\right)  }{m+n+1}%
+\dfrac{2\cdot\left(  2n+1\right)  }{m+n+1}\right)  }_{\substack{=4\\\text{(by
straightforward computation)}}}\cdot T\left(  m,n\right)  =4T\left(
m,n\right)  .
\end{align*}
In other words, $4T\left(  m,n\right)  =T\left(  m+1,n\right)  +T\left(
m,n+1\right)  $. This solves Exercise \ref{exe.supercat} \textbf{(a)}.
\end{vershort}

\begin{verlong}
\textbf{(a)} The definition of $T$ yields%
\begin{equation}
T\left(  m,n\right)  =\dfrac{\left(  2m\right)  !\left(  2n\right)
!}{m!n!\left(  m+n\right)  !}. \label{sol.supercat.a.0}%
\end{equation}

Every $q\in\mathbb{N}$ satisfies%
\begin{equation}
\left(  q+1\right)  !=q!\cdot\left(  q+1\right)  . \label{sol.supercat.a.1}%
\end{equation}

Let $m\in\mathbb{N}$ and $n\in\mathbb{N}$. Then, $2\left(  m+1\right)
=2m+2=\left(  2m+1\right)  +1$, so that%
\begin{align}
\left(  2\left(  m+1\right)  \right)  !  &  =\left(  \left(  2m+1\right)
+1\right)  !=\underbrace{\left(  2m+1\right)  !}_{\substack{=\left(
2m\right)  !\cdot\left(  2m+1\right)  \\\text{(by (\ref{sol.supercat.a.1})
(applied to }q=2m\text{))}}}\cdot\underbrace{\left(  \left(  2m+1\right)
+1\right)  }_{=2\left(  m+1\right)  }\nonumber\\
&  \ \ \ \ \ \ \ \ \ \ \left(  \text{by (\ref{sol.supercat.a.1}) (applied to
}q=2m+1\text{)}\right) \nonumber\\
&  =\left(  2m\right)  !\cdot\left(  2m+1\right)  \cdot2\left(  m+1\right)
\nonumber\\
&  =2\cdot\left(  2m\right)  !\cdot\left(  2m+1\right)  \left(  m+1\right)  .
\label{sol.supercat.a.2}%
\end{align}

Furthermore,%
\begin{equation}
\left(  m+1\right)  !=m!\cdot\left(  m+1\right)  \label{sol.supercat.a.3}%
\end{equation}
(by (\ref{sol.supercat.a.1}) (applied to $q=m$)). Moreover, $\left(
m+1\right)  +n=\left(  m+n\right)  +1$, so that
\begin{equation}
\left(  \left(  m+1\right)  +n\right)  !=\left(  \left(  m+n\right)
+1\right)  !=\left(  m+n\right)  !\cdot\left(  \left(  m+n\right)  +1\right)
\label{sol.supercat.a.4}%
\end{equation}
(by (\ref{sol.supercat.a.1}) (applied to $q=m+n$)).

On the other hand, the definition of $T\left(  m+1,n\right)  $ yields%
\[
T\left(  m+1,n\right)  =\dfrac{\left(  2\left(  m+1\right)  \right)  !\left(
2n\right)  !}{\left(  m+1\right)  !n!\left(  \left(  m+1\right)  +n\right)
!}.
\]
In view of (\ref{sol.supercat.a.2}), (\ref{sol.supercat.a.3}) and
(\ref{sol.supercat.a.4}), this rewrites as%
\[
T\left(  m+1,n\right)  =\dfrac{2\cdot\left(  2m\right)  !\cdot\left(
2m+1\right)  \left(  m+1\right)  \cdot\left(  2n\right)  !}{m!\cdot\left(
m+1\right)  \cdot n!\cdot\left(  m+n\right)  !\cdot\left(  \left(  m+n\right)
+1\right)  }.
\]
Hence,%
\begin{align}
T\left(  m+1,n\right)   &  =\dfrac{2\cdot\left(  2m\right)  !\cdot\left(
2m+1\right)  \left(  m+1\right)  \cdot\left(  2n\right)  !}{m!\cdot\left(
m+1\right)  \cdot n!\cdot\left(  m+n\right)  !\cdot\left(  \left(  m+n\right)
+1\right)  }\nonumber\\
&  =\dfrac{2\cdot\left(  2m+1\right)  }{\left(  m+n\right)  +1}%
\underbrace{\dfrac{\left(  2m\right)  !\left(  2n\right)  !}{m!n!\left(
m+n\right)  !}}_{\substack{=T\left(  m,n\right)  \\\text{(by
(\ref{sol.supercat.a.0}))}}}\nonumber\\
&  =\dfrac{2\cdot\left(  2m+1\right)  }{\left(  m+n\right)  +1}T\left(
m,n\right)  . \label{sol.supercat.a.5}%
\end{align}
The same argument (but with $m$ and $n$ replaced by $n$ and $m$) shows that%
\begin{equation}
T\left(  n+1,m\right)  =\dfrac{2\cdot\left(  2n+1\right)  }{\left(
n+m\right)  +1}T\left(  n,m\right)  . \label{sol.supercat.a.6}%
\end{equation}

But Exercise \ref{exe.supercat} \textbf{(g)} yields $T\left(  m,n\right)
=T\left(  n,m\right)  $. Hence, $T\left(  n,m\right)  =T\left(  m,n\right)  $.
Also, Exercise \ref{exe.supercat} \textbf{(g)} (applied to $n+1$ instead of
$n$) yields%
\begin{align*}
T\left(  m,n+1\right)   &  =T\left(  n+1,m\right)  =\dfrac{2\cdot\left(
2n+1\right)  }{\left(  n+m\right)  +1}\underbrace{T\left(  n,m\right)
}_{=T\left(  m,n\right)  }=\dfrac{2\cdot\left(  2n+1\right)  }{\left(
n+m\right)  +1}T\left(  m,n\right) \\
&  =\dfrac{2\cdot\left(  2n+1\right)  }{\left(  m+n\right)  +1}T\left(
m,n\right)
\end{align*}
(since $n+m=m+n$). Adding this equality to (\ref{sol.supercat.a.5}), we obtain%
\begin{align*}
T\left(  m+1,n\right)  +T\left(  m,n+1\right)   &  =\dfrac{2\cdot\left(
2m+1\right)  }{\left(  m+n\right)  +1}T\left(  m,n\right)  +\dfrac
{2\cdot\left(  2n+1\right)  }{\left(  m+n\right)  +1}T\left(  m,n\right) \\
&  =\underbrace{\left(  \dfrac{2\cdot\left(  2m+1\right)  }{\left(
m+n\right)  +1}+\dfrac{2\cdot\left(  2n+1\right)  }{\left(  m+n\right)
+1}\right)  }_{\substack{=4\\\text{(by straightforward computation)}}}T\left(
m,n\right) \\
&  =4T\left(  m,n\right)  .
\end{align*}
In other words, $4T\left(  m,n\right)  =T\left(  m+1,n\right)  +T\left(
m,n+1\right)  $. This solves Exercise \ref{exe.supercat} \textbf{(a)}.
\end{verlong}

Before we come to the solution of Exercise \ref{exe.supercat} \textbf{(b)},
let us observe something trivial:

\begin{statement}
\textit{Observation 3:} Let $m\in\mathbb{N}$ and $n\in\mathbb{N}$. Then,
$T\left(  m,n\right)  >0$.
\end{statement}

\begin{vershort}
[\textit{Proof of Observation 3:} This follows immediately from the definition
of $T\left(  m,n\right)  $ (since the numbers $\left(  2m\right)  !$, $\left(
2n\right)  !$, $m!$, $n!$ and $\left(  m+n\right)  !$ are positive).]
\end{vershort}

\begin{verlong}
[\textit{Proof of Observation 3:} The definition of $T\left(  m,n\right)  $
yields $T\left(  m,n\right)  =\dfrac{\left(  2m\right)  !\left(  2n\right)
!}{m!n!\left(  m+n\right)  !}>0$ (since $\left(  2m\right)  !$, $\left(
2n\right)  !$, $m!$, $n!$ and $\left(  m+n\right)  !$ are positive integers
(because $q!$ is a positive integer for every $q\in\mathbb{N}$)). This proves
Observation 3.]
\end{verlong}

\begin{vershort}
\textbf{(b)} \textit{First solution to Exercise \ref{exe.supercat}
\textbf{(b)}:} Let $m\in\mathbb{N}$ and $n\in\mathbb{N}$. Let $p=\min\left\{
m,n\right\}  $. Exercise \ref{exe.supercat} \textbf{(i)} yields%
\[
T\left(  m,n\right)  =\sum_{k=-p}^{p}\left(  -1\right)  ^{k}\dbinom{2m}%
{m+k}\dbinom{2n}{n-k}.
\]
The right-hand side of this equality is clearly an integer (since the binomial
coefficients appearing in it are integers\footnote{according to Proposition
\ref{prop.binom.int}}). Thus, so is the left-hand side. In other words,
$T\left(  m,n\right)  $ is an integer. But Observation 3 yields $T\left(
m,n\right)  >0$. Hence, $T\left(  m,n\right)  $ is a positive integer (since
$T\left(  m,n\right)  $ is an integer). Therefore, $T\left(  m,n\right)
\in\mathbb{N}$. This solves Exercise \ref{exe.supercat} \textbf{(b)}.
\end{vershort}

\begin{verlong}
\textbf{(b)} \textit{First solution to Exercise \ref{exe.supercat}
\textbf{(b)}:} Let $m\in\mathbb{N}$ and $n\in\mathbb{N}$. Let $p=\min\left\{
m,n\right\}  $. Exercise \ref{exe.supercat} \textbf{(i)} yields%
\begin{equation}
T\left(  m,n\right)  =\sum_{k=-p}^{p}\left(  -1\right)  ^{k}\dbinom{2m}%
{m+k}\dbinom{2n}{n-k}. \label{sol.supercat.b.1}%
\end{equation}

For every $k\in\left\{  -p,-p+1,\ldots,p\right\}  $, the number $\left(
-1\right)  ^{k}\dbinom{2m}{m+k}\dbinom{2n}{n-k}$ is an
integer\footnote{\textit{Proof.} Let $k\in\left\{  -p,-p+1,\ldots,p\right\}
$. We must prove that $\left(  -1\right)  ^{k}\dbinom{2m}{m+k}\dbinom{2n}%
{n-k}$ is an integer.
\par
Clearly, $p=\min\left\{  m,n\right\}  \in\mathbb{N}$ (since $m\in\mathbb{N}$
and $n\in\mathbb{N}$). We have $k\in\left\{  -p,-p+1,\ldots,p\right\}  $.
Thus, $k$ is an integer satisfying $-p\leq k\leq p$.
\par
We have $p=\min\left\{  m,n\right\}  \leq m$. Also, $-p\leq k$, so that
$k\geq-\underbrace{p}_{\leq m}\geq-m$. Hence, $m+k\geq0$. Hence,
$m+k\in\mathbb{N}$ (since $m+k$ is an integer (since $m$ and $k$ are
integers)).
\par
We have $p=\min\left\{  m,n\right\}  \leq n$. Hence, $k\leq p\leq n$. Thus,
$n-k\geq0$. Hence, $n-k\in\mathbb{N}$ (since $n-k$ is an integer (since $n$
and $k$ are integers)).
\par
Also, $m\in\mathbb{N}$ and thus $2m\in\mathbb{N}\subseteq\mathbb{Z}$. Hence,
Proposition \ref{prop.binom.int} (applied to $2m$ and $m+k$ instead of $m$ and
$n$) yields $\dbinom{2m}{m+k}\in\mathbb{Z}$. In other words, $\dbinom{2m}%
{m+k}$ is an integer.
\par
Furthermore, $n\in\mathbb{N}$ and thus $2n\in\mathbb{N}\subseteq\mathbb{Z}$.
Hence, Proposition \ref{prop.binom.int} (applied to $2n$ and $n-k$ instead of
$m$ and $n$) yields $\dbinom{2n}{n-k}\in\mathbb{Z}$. In other words,
$\dbinom{2n}{n-k}$ is an integer.
\par
Finally, $\left(  -1\right)  ^{k}$ is an integer (since $k$ is an integer).
\par
Now, we know that all three numbers $\left(  -1\right)  ^{k}$, $\dbinom
{2m}{m+k}$ and $\dbinom{2n}{n-k}$ are integers. Hence, their product $\left(
-1\right)  ^{k}\dbinom{2m}{m+k}\dbinom{2n}{n-k}$ is an integer as well. Qed.}.
Hence, the sum $\sum_{k=-p}^{p}\left(  -1\right)  ^{k}\dbinom{2m}{m+k}%
\dbinom{2n}{n-k}$ is a sum of integers, and thus itself an integer. In other
words, $\sum_{k=-p}^{p}\left(  -1\right)  ^{k}\dbinom{2m}{m+k}\dbinom{2n}%
{n-k}\in\mathbb{Z}$. In view of (\ref{sol.supercat.b.1}), this rewrites as
$T\left(  m,n\right)  \in\mathbb{Z}$. In other words, $T\left(  m,n\right)  $
is an integer.

But Observation 3 yields $T\left(  m,n\right)  >0$. In other words, $T\left(
m,n\right)  $ is positive. Hence, $T\left(  m,n\right)  $ is a positive
integer. Therefore, $T\left(  m,n\right)  \in\mathbb{N}$. This solves Exercise
\ref{exe.supercat} \textbf{(b)}.
\end{verlong}

\textit{Second solution to Exercise \ref{exe.supercat} \textbf{(b)}:} We shall
prove Exercise \ref{exe.supercat} \textbf{(b)} by induction on $n$:

\begin{vershort}
\textit{Induction base:} We have $T\left(  m,0\right)  \in\mathbb{N}$ for
every $m\in\mathbb{N}$\ \ \ \ \footnote{\textit{Proof.} Let $m\in\mathbb{N}$.
Then, Exercise \ref{exe.supercat} \textbf{(e)} shows that $T\left(
m,0\right)  =\dbinom{2m}{m}\in\mathbb{N}$ (this follows from an application of
Lemma \ref{lem.binom.intN}). Qed.}. In other words, Exercise
\ref{exe.supercat} \textbf{(b)} holds for $n=0$. This completes the induction base.
\end{vershort}

\begin{verlong}
\textit{Induction base:} We have $T\left(  m,0\right)  \in\mathbb{N}$ for
every $m\in\mathbb{N}$\ \ \ \ \footnote{\textit{Proof.} Let $m\in\mathbb{N}$.
Then, Proposition \ref{prop.binom.int} (applied to $2m$ and $m$ instead of $m$
and $n$) yields $\dbinom{2m}{m}\in\mathbb{Z}$. Now, Exercise
\ref{exe.supercat} \textbf{(e)} shows that $T\left(  m,0\right)  =\dbinom
{2m}{m}\in\mathbb{Z}$. But Observation 3 (applied to $n=0$) yields $T\left(
m,0\right)  >0$. Thus, $T\left(  m,0\right)  $ is positive. Combining this
with $T\left(  m,0\right)  \in\mathbb{Z}$, we conclude that $T\left(
m,0\right)  $ is a positive integer. Hence, $T\left(  m,0\right)
\in\mathbb{N}$. Qed.}. In other words, Exercise \ref{exe.supercat}
\textbf{(b)} holds for $n=0$. This completes the induction base.
\end{verlong}

\textit{Induction step:} Let $N\in\mathbb{N}$. Assume that Exercise
\ref{exe.supercat} \textbf{(b)} holds for $n=N$. We must prove that Exercise
\ref{exe.supercat} \textbf{(b)} holds for $n=N+1$.

We have assumed that Exercise \ref{exe.supercat} \textbf{(b)} holds for $n=N$.
In other words, we have
\begin{equation}
T\left(  m,N\right)  \in\mathbb{N}\ \ \ \ \ \ \ \ \ \ \text{for every }%
m\in\mathbb{N}. \label{sol.supercat.b.indhyp}%
\end{equation}

Now, let $m\in\mathbb{N}$. We shall show that $T\left(  m,N+1\right)
\in\mathbb{N}$.

Indeed, (\ref{sol.supercat.b.indhyp}) yields $T\left(  m,N\right)
\in\mathbb{N}\subseteq\mathbb{Z}$. But (\ref{sol.supercat.b.indhyp}) (applied
to $m+1$ instead of $m$) yields $T\left(  m+1,N\right)  \in\mathbb{N}%
\subseteq\mathbb{Z}$. Hence, $4T\left(  m,N\right)  -T\left(  m+1,N\right)
\in\mathbb{Z}$ (since both $T\left(  m,N\right)  $ and $T\left(  m+1,N\right)
$ belong to $\mathbb{Z}$).

Exercise \ref{exe.supercat} \textbf{(a)} (applied to $n=N$) yields $4T\left(
m,N\right)  =T\left(  m+1,N\right)  +T\left(  m,N+1\right)  $. Hence,
$T\left(  m,N+1\right)  =4T\left(  m,N\right)  -T\left(  m+1,N\right)
\in\mathbb{Z}$.

\begin{vershort}
But Observation 3 (applied to $n=N+1$) yields $T\left(  m,N+1\right)  >0$.
Combining this with $T\left(  m,N+1\right)  \in\mathbb{Z}$, we conclude that
$T\left(  m,N+1\right)  \in\mathbb{N}$.
\end{vershort}

\begin{verlong}
Hence, $T\left(  m,N+1\right)  $ is an integer. But Observation 3 (applied to
$n=N+1$) yields $T\left(  m,N+1\right)  >0$. Thus, $T\left(  m,N+1\right)  $
is a positive integer (since $T\left(  m,N+1\right)  $ is an integer).
Therefore, $T\left(  m,N+1\right)  \in\mathbb{N}$.
\end{verlong}

Now, forget that we fixed $m$. We thus have proven that $T\left(
m,N+1\right)  \in\mathbb{N}$ for every $m\in\mathbb{N}$. In other words,
Exercise \ref{exe.supercat} \textbf{(b)} holds for $n=N+1$. This completes the
induction step. Thus, the induction proof of Exercise \ref{exe.supercat}
\textbf{(b)} is complete.

\textbf{(c)} For every $m\in\mathbb{N}$ and $n\in\mathbb{N}$, we have
$T\left(  m,n\right)  \in\mathbb{N}$ (by Exercise \ref{exe.supercat}
\textbf{(b)}). Thus, for every $m\in\mathbb{N}$ and $n\in\mathbb{N}$, the
number $T\left(  m,n\right)  $ is an integer. Hence, speaking of
\textquotedblleft the integer $T\left(  m,n\right)  $\textquotedblright\ in
Exercise \ref{exe.supercat} \textbf{(c)} makes sense.

We shall prove Exercise \ref{exe.supercat} \textbf{(c)} by induction on $n$:

\textit{Induction base:} If $m\in\mathbb{N}$ is such that $\left(  m,0\right)
\neq\left(  0,0\right)  $, then the integer $T\left(  m,0\right)  $ is
even\footnote{\textit{Proof.} Let $m\in\mathbb{N}$ be such that $\left(
m,0\right)  \neq\left(  0,0\right)  $. Then, Exercise \ref{exe.supercat}
\textbf{(e)} shows that $T\left(  m,0\right)  =\dbinom{2m}{m}$.
\par
If we had $m=0$, then we would have $\left(  \underbrace{m}_{=0},0\right)
=\left(  0,0\right)  $, which would contradict $\left(  m,0\right)
\neq\left(  0,0\right)  $. Hence, we cannot have $m=0$. Thus, we must have
$m\neq0$. Therefore, $m$ is a positive integer (since $m\in\mathbb{N}$).
Hence, Exercise \ref{exe.central-binomial-even} \textbf{(a)} shows that the
binomial coefficient $\dbinom{2m}{m}$ is even. In other words, the integer
$\dbinom{2m}{m}$ is even. In view of $T\left(  m,0\right)  =\dbinom{2m}{m}$,
this rewrites as follows: The integer $T\left(  m,0\right)  $ is even. Qed.}.
In other words, Exercise \ref{exe.supercat} \textbf{(c)} holds for $n=0$. This
completes the induction base.

\textit{Induction step:} Let $N\in\mathbb{N}$. Assume that Exercise
\ref{exe.supercat} \textbf{(c)} holds for $n=N$. We must prove that Exercise
\ref{exe.supercat} \textbf{(c)} holds for $n=N+1$.

We have assumed that Exercise \ref{exe.supercat} \textbf{(c)} holds for $n=N$.
In other words, if $m\in\mathbb{N}$ is such that $\left(  m,N\right)
\neq\left(  0,0\right)  $, then%
\begin{equation}
\text{the integer }T\left(  m,N\right)  \text{ is even.}
\label{sol.supercat.c.indhyp}%
\end{equation}

Now, let $m\in\mathbb{N}$ be such that $\left(  m,N+1\right)  \neq\left(
0,0\right)  $. We shall show that the integer $T\left(  m,N+1\right)  $ is even.

\begin{vershort}
We have $\left(  m+1,N\right)  \neq\left(  0,0\right)  $ (since $m+1\neq0$).
Hence, (\ref{sol.supercat.c.indhyp}) (applied to $m+1$ instead of $m$) shows
that the integer $T\left(  m+1,N\right)  $ is even. In other words, $T\left(
m+1,N\right)  \equiv0\operatorname{mod}2$.
\end{vershort}

\begin{verlong}
We have $\left(  m+1,N\right)  \neq\left(  0,0\right)  $%
\ \ \ \ \footnote{\textit{Proof.} Assume the contrary. Thus, $\left(
m+1,N\right)  =\left(  0,0\right)  $. Hence, $m+1=0$ and $N=0$. In particular,
$m+1=0$, so that $m=-1\notin\mathbb{N}$. This contradicts $m\in\mathbb{N}$.
\par
This contradiction shows that our assumption was wrong. Qed.}. Hence,
(\ref{sol.supercat.c.indhyp}) (applied to $m+1$ instead of $m$) shows that the
integer $T\left(  m+1,N\right)  $ is even. In other words, $T\left(
m+1,N\right)  \equiv0\operatorname{mod}2$.
\end{verlong}

\begin{vershort}
But Exercise \ref{exe.supercat} \textbf{(b)} (applied to $n=N$) yields
$T\left(  m,N\right)  \in\mathbb{N}$. Hence, $4T\left(  m,N\right)
\equiv0\operatorname{mod}2$.
\end{vershort}

\begin{verlong}
But Exercise \ref{exe.supercat} \textbf{(b)} (applied to $n=N$) yields
$T\left(  m,N\right)  \in\mathbb{N}$. Hence, $2T\left(  m,N\right)
\equiv0\operatorname{mod}2$. Thus, $4T\left(  m,N\right)  =2\cdot
\underbrace{2T\left(  m,N\right)  }_{\equiv0\operatorname{mod}2}%
\equiv0\operatorname{mod}2$.
\end{verlong}

Exercise \ref{exe.supercat} \textbf{(a)} (applied to $n=N$) yields $4T\left(
m,N\right)  =T\left(  m+1,N\right)  +T\left(  m,N+1\right)  $. Hence,
\[
T\left(  m,N+1\right)  =\underbrace{4T\left(  m,N\right)  }_{\equiv
0\operatorname{mod}2}-\underbrace{T\left(  m+1,N\right)  }_{\equiv
0\operatorname{mod}2}\equiv0-0=0\operatorname{mod}2.
\]
In other words, the integer $T\left(  m,N+1\right)  $ is even.

Now, forget that we fixed $m$. We thus have proven that if $m\in\mathbb{N}$ is
such that $\left(  m,N+1\right)  \neq\left(  0,0\right)  $, then the integer
$T\left(  m,N+1\right)  $ is even. In other words, Exercise \ref{exe.supercat}
\textbf{(c)} holds for $n=N+1$. This completes the induction step. Thus, the
induction proof of Exercise \ref{exe.supercat} \textbf{(c)} is complete.

\textbf{(d)} For every $m\in\mathbb{N}$ and $n\in\mathbb{N}$, we have
$T\left(  m,n\right)  \in\mathbb{N}$ (by Exercise \ref{exe.supercat}
\textbf{(b)}). Thus, for every $m\in\mathbb{N}$ and $n\in\mathbb{N}$, the
number $T\left(  m,n\right)  $ is an integer. Hence, the divisibility
statement \textquotedblleft$4\mid T\left(  m,n\right)  $\textquotedblright\ in
Exercise \ref{exe.supercat} \textbf{(d)} makes sense.

We shall prove Exercise \ref{exe.supercat} \textbf{(d)} by induction on $n$:

\textit{Induction base:} If $m\in\mathbb{N}$ is such that $m+0$ is odd and
$m+0>1$, then $4\mid T\left(  m,0\right)  $\ \ \ \ \footnote{\textit{Proof.}
Let $m\in\mathbb{N}$ be such that $m+0$ is odd and $m+0>1$. Then, Exercise
\ref{exe.supercat} \textbf{(e)} shows that $T\left(  m,0\right)  =\dbinom
{2m}{m}$.
\par
We know that $m+0$ is odd. In other words, $m$ is odd (since $m=m+0$).
Furthermore, $m$ is a positive integer (since $m=m+0>1>0$). Therefore,
Exercise \ref{exe.central-binomial-even} \textbf{(c)} shows that $\dbinom
{2m}{m}\equiv0\operatorname{mod}4$. Thus, $T\left(  m,0\right)  =\dbinom
{2m}{m}\equiv0\operatorname{mod}4$. In other words, $4\mid T\left(
m,0\right)  $. Qed.}. In other words, Exercise \ref{exe.supercat} \textbf{(d)}
holds for $n=0$. This completes the induction base.

\textit{Induction step:} Let $N\in\mathbb{N}$. Assume that Exercise
\ref{exe.supercat} \textbf{(d)} holds for $n=N$. We must prove that Exercise
\ref{exe.supercat} \textbf{(d)} holds for $n=N+1$.

We have assumed that Exercise \ref{exe.supercat} \textbf{(d)} holds for $n=N$.
In other words, if $m\in\mathbb{N}$ is such that $m+N$ is odd and $m+N>1$,
then%
\begin{equation}
4\mid T\left(  m,N\right)  \text{.} \label{sol.supercat.d.indhyp}%
\end{equation}

Now, let $m\in\mathbb{N}$ be such that $m+\left(  N+1\right)  $ is odd and
$m+\left(  N+1\right)  >1$. We shall show that $4\mid T\left(  m,N+1\right)  $.

We know that $m+\left(  N+1\right)  $ is odd. In other words, $\left(
m+1\right)  +N$ is odd (since $\left(  m+1\right)  +N=m+\left(  N+1\right)
$). Also, $\left(  m+1\right)  +N=m+\left(  N+1\right)  >1$. Hence,
(\ref{sol.supercat.d.indhyp}) (applied to $m+1$ instead of $m$) shows that
$4\mid T\left(  m+1,N\right)  $. In other words, $T\left(  m+1,N\right)
\equiv0\operatorname{mod}4$.

But Exercise \ref{exe.supercat} \textbf{(b)} (applied to $n=N$) yields
$T\left(  m,N\right)  \in\mathbb{N}$. Hence, $4T\left(  m,N\right)
\equiv0\operatorname{mod}4$.

Exercise \ref{exe.supercat} \textbf{(a)} (applied to $n=N$) yields $4T\left(
m,N\right)  =T\left(  m+1,N\right)  +T\left(  m,N+1\right)  $. Hence,
\[
T\left(  m,N+1\right)  =\underbrace{4T\left(  m,N\right)  }_{\equiv
0\operatorname{mod}4}-\underbrace{T\left(  m+1,N\right)  }_{\equiv
0\operatorname{mod}4}\equiv0-0=0\operatorname{mod}4.
\]
In other words, $4\mid T\left(  m,N+1\right)  $.

Now, forget that we fixed $m$. We thus have proven that if $m\in\mathbb{N}$ is
such that $m+\left(  N+1\right)  $ is odd and $m+\left(  N+1\right)  >1$, then
$4\mid T\left(  m,N+1\right)  $. In other words, Exercise \ref{exe.supercat}
\textbf{(d)} holds for $n=N+1$. This completes the induction step. Thus, the
induction proof of Exercise \ref{exe.supercat} \textbf{(d)} is complete.
\end{proof}

\subsection{Solution to Exercise \ref{exe.choose.a/b}}

\subsubsection{First solution}

Before we solve Exercise \ref{exe.choose.a/b}, let us prove two basic facts in
modular arithmetic:

\begin{proposition}
\label{prop.sol.choose.a/b.lem1}Let $b$ and $c$ be two integers such that
$c>0$. Then, there exists an $s\in\mathbb{Z}$ such that $b^{c-1}\equiv
sb^{c}\operatorname{mod}c$.
\end{proposition}

\begin{proof}
[Proof of Proposition \ref{prop.sol.choose.a/b.lem1}.]For every integer $m$,
let $m\%c$ denote the remainder obtained when $m$ is divided by $c$. Thus,
every integer $m$ satisfies%
\begin{equation}
m\%c\in\left\{  0,1,\ldots,c-1\right\}
\label{pf.prop.sol.choose.a/b.lem1.rem.1}%
\end{equation}
and%
\begin{equation}
m\%c\equiv m\operatorname{mod}c. \label{pf.prop.sol.choose.a/b.lem1.rem.2}%
\end{equation}
Indeed, these two relations follow from Corollary \ref{cor.ind.quo-rem.remmod}
\textbf{(a)} (applied to $N=c$ and $n=m$).

We shall now show that two of the $c+1$ integers $b^{0}\%c,b^{1}%
\%c,\ldots,b^{c}\%c$ are equal.

Indeed, assume the contrary (for the sake of contradiction). Thus, the $c+1$
integers $b^{0}\%c,b^{1}\%c,\ldots,b^{c}\%c$ are pairwise distinct. Hence,%
\[
\left\vert \left\{  b^{0}\%c,b^{1}\%c,\ldots,b^{c}\%c\right\}  \right\vert
=c+1.
\]
But we have $b^{i}\%c\in\left\{  0,1,\ldots,c-1\right\}  $ for each
$i\in\left\{  0,1,\ldots,c\right\}  $ (by
(\ref{pf.prop.sol.choose.a/b.lem1.rem.1}), applied to $m=b^{i}$). Thus,%
\[
\left\{  b^{0}\%c,b^{1}\%c,\ldots,b^{c}\%c\right\}  \subseteq\left\{
0,1,\ldots,c-1\right\}  .
\]
Therefore,%
\[
\left\vert \left\{  b^{0}\%c,b^{1}\%c,\ldots,b^{c}\%c\right\}  \right\vert
\leq\left\vert \left\{  0,1,\ldots,c-1\right\}  \right\vert =c,
\]
so that%
\[
c\geq\left\vert \left\{  b^{0}\%c,b^{1}\%c,\ldots,b^{c}\%c\right\}
\right\vert =c+1.
\]
This contradicts $c<c+1$. This contradiction shows that our assumption was
wrong. Hence, we have proven that two of the $c+1$ integers $b^{0}%
\%c,b^{1}\%c,\ldots,b^{c}\%c$ are equal. In other words, there exist two
distinct elements $u$ and $v$ of $\left\{  0,1,\ldots,c\right\}  $ such that
$b^{u}\%c=b^{v}\%c$. Consider these $u$ and $v$.

We can WLOG assume that $u\leq v$ (since otherwise, we can simply swap $u$
with $v$, and nothing changes). Assume this. Thus, $u\leq v$, so that $u<v$
(since $u$ and $v$ are distinct).

Now, (\ref{pf.prop.sol.choose.a/b.lem1.rem.2}) (applied to $m=b^{u}$) yields
$b^{u}\%c\equiv b^{u}\operatorname{mod}c$. Also,
(\ref{pf.prop.sol.choose.a/b.lem1.rem.2}) (applied to $m=b^{v}$) yields
$b^{v}\%c\equiv b^{v}\operatorname{mod}c$. Thus, $b^{u}\equiv b^{u}%
\%c=b^{v}\%c\equiv b^{v}\operatorname{mod}c$.

Also, $u\in\left\{  0,1,\ldots,c\right\}  $, so that $0\leq u$. Also,
$v\in\left\{  0,1,\ldots,c\right\}  $, so that $v\leq c$. Hence,
$c-v\in\mathbb{N}$. Thus, $b^{c-v}$ is a well-defined integer.

Also, $0\leq u<v$, so that $0\leq v-1$ (since $0$ and $v$ are integers). Thus,
$v-1\in\mathbb{N}$, so that $b^{v-1}$ is a well-defined integer.

But $u<v$, and thus $u\leq v-1$ (since $u$ and $v$ are integers). Hence,
$\left(  v-1\right)  -u\in\mathbb{N}$. Thus, $b^{\left(  v-1\right)  -u}$ is a
well-defined integer. Set $t=b^{\left(  v-1\right)  -u}$. Thus, $t$ is an integer.

Now, $v-1=\left(  \left(  v-1\right)  -u\right)  +u$, and thus%
\begin{align*}
b^{v-1}  &  =b^{\left(  \left(  v-1\right)  -u\right)  +u}%
=\underbrace{b^{\left(  v-1\right)  -u}}_{=t}b^{u}\ \ \ \ \ \ \ \ \ \ \left(
\text{since }0\leq u\leq v-1\right) \\
&  =t\underbrace{b^{u}}_{\equiv b^{v}\operatorname{mod}c}\equiv tb^{v}%
\operatorname{mod}c.
\end{align*}

Now, $c-1=\left(  v-1\right)  +\left(  c-v\right)  $, so that%
\begin{align*}
b^{c-1}  &  =b^{\left(  v-1\right)  +\left(  c-v\right)  }=\underbrace{b^{v-1}%
}_{\equiv tb^{v}\operatorname{mod}c}b^{c-v}\ \ \ \ \ \ \ \ \ \ \left(
\text{since }v-1\in\mathbb{N}\text{ and }c-v\in\mathbb{N}\right) \\
&  \equiv t\underbrace{b^{v}b^{c-v}}_{=b^{v+\left(  c-v\right)  }=b^{c}%
}=tb^{c}\operatorname{mod}c.
\end{align*}
Hence, there exists an $s\in\mathbb{Z}$ such that $b^{c-1}\equiv
sb^{c}\operatorname{mod}c$ (namely, $s=t$). This proves Proposition
\ref{prop.sol.choose.a/b.lem1}.
\end{proof}

\begin{lemma}
\label{lem.sol.choose.a/b.prod-of-congs}Let $n\in\mathbb{N}$. Let $c$ be a
positive integer. Let $u_{0},u_{1},\ldots,u_{n-1}$ be $n$ integers. Let
$v_{0},v_{1},\ldots,v_{n-1}$ be $n$ integers. Let $d$ be an integer. Assume
that%
\begin{equation}
du_{i}\equiv dv_{i}\operatorname{mod}c\ \ \ \ \ \ \ \ \ \ \text{for each }%
i\in\left\{  0,1,\ldots,n-1\right\}  .
\label{eq.lem.sol.choose.a/b.prod-of-congs.cond}%
\end{equation}
Then,%
\[
d\prod_{i=0}^{n-1}u_{i}\equiv d\prod_{i=0}^{n-1}v_{i}\operatorname{mod}c.
\]

\end{lemma}

\begin{proof}
[Proof of Lemma \ref{lem.sol.choose.a/b.prod-of-congs}.]We claim that%
\begin{equation}
d\prod_{i=0}^{k-1}u_{i}\equiv d\prod_{i=0}^{k-1}v_{i}\operatorname{mod}c
\label{pf.lem.sol.choose.a/b.prod-of-congs.goal}%
\end{equation}
for each $k\in\left\{  0,1,\ldots,n\right\}  $.

[\textit{Proof of (\ref{pf.lem.sol.choose.a/b.prod-of-congs.goal}):} We shall
prove (\ref{pf.lem.sol.choose.a/b.prod-of-congs.goal}) by induction over $k$:

\textit{Induction base:} We have $d\underbrace{\prod_{i=0}^{0-1}u_{i}%
}_{=\left(  \text{empty product}\right)  =1}=d$. Comparing this with
\newline$d\underbrace{\prod_{i=0}^{0-1}v_{i}}_{=\left(  \text{empty
product}\right)  =1}=d$, we obtain $d\prod_{i=0}^{0-1}u_{i}=d\prod_{i=0}%
^{0-1}v_{i}$. Hence, $d\prod_{i=0}^{0-1}u_{i}\equiv d\prod_{i=0}^{0-1}%
v_{i}\operatorname{mod}c$. In other words,
(\ref{pf.lem.sol.choose.a/b.prod-of-congs.goal}) holds for $k=0$. This
completes the induction base.

\textit{Induction step:} Let $K\in\left\{  0,1,\ldots,n\right\}  $ be
positive. Assume that (\ref{pf.lem.sol.choose.a/b.prod-of-congs.goal}) holds
for $k=K-1$. We must show that (\ref{pf.lem.sol.choose.a/b.prod-of-congs.goal}%
) holds for $k=K$.

We have $K\in\left\{  1,2,\ldots,n\right\}  $ (since $K\in\left\{
0,1,\ldots,n\right\}  $ and $K$ is positive). Hence, $K-1\in\left\{
0,1,\ldots,n-1\right\}  $.

We have assumed that (\ref{pf.lem.sol.choose.a/b.prod-of-congs.goal}) holds
for $k=K-1$. In other words, we have%
\begin{equation}
d\prod_{i=0}^{\left(  K-1\right)  -1}u_{i}\equiv d\prod_{i=0}^{\left(
K-1\right)  -1}v_{i}\operatorname{mod}c.
\label{pf.lem.sol.choose.a/b.prod-of-congs.goal.pf.1}%
\end{equation}

Now,%
\begin{align*}
d\underbrace{\prod_{i=0}^{K-1}u_{i}}_{=\left(  \prod_{i=0}^{\left(
K-1\right)  -1}u_{i}\right)  u_{K-1}}  &  =\underbrace{d\left(  \prod
_{i=0}^{\left(  K-1\right)  -1}u_{i}\right)  }_{\substack{\equiv d\prod
_{i=0}^{\left(  K-1\right)  -1}v_{i}\operatorname{mod}c\\\text{(by
(\ref{pf.lem.sol.choose.a/b.prod-of-congs.goal.pf.1}))}}}u_{K-1}\\
&  \equiv\left(  d\prod_{i=0}^{\left(  K-1\right)  -1}v_{i}\right)
u_{K-1}=\underbrace{du_{K-1}}_{\substack{\equiv dv_{K-1}\operatorname{mod}%
c\\\text{(by (\ref{eq.lem.sol.choose.a/b.prod-of-congs.cond}), applied to
}i=K-1\text{)}}}\prod_{i=0}^{\left(  K-1\right)  -1}v_{i}\\
&  \equiv dv_{K-1}\prod_{i=0}^{\left(  K-1\right)  -1}v_{i}%
=d\underbrace{\left(  \prod_{i=0}^{\left(  K-1\right)  -1}v_{i}\right)
v_{K-1}}_{=\prod_{i=0}^{K-1}v_{i}}=d\prod_{i=0}^{K-1}v_{i}\operatorname{mod}c.
\end{align*}
In other words, (\ref{pf.lem.sol.choose.a/b.prod-of-congs.goal}) holds for
$k=K$. This completes the induction step. Thus,
(\ref{pf.lem.sol.choose.a/b.prod-of-congs.goal}) is proven by induction.]

Now, (\ref{pf.lem.sol.choose.a/b.prod-of-congs.goal}) (applied to $k=n$)
yields%
\[
d\prod_{i=0}^{n-1}u_{i}\equiv d\prod_{i=0}^{n-1}v_{i}\operatorname{mod}c.
\]
This proves Lemma \ref{lem.sol.choose.a/b.prod-of-congs}.
\end{proof}

\begin{proof}
[First solution to Exercise \ref{exe.choose.a/b}.]Set $c=n!$. Notice that $c$
is a positive integer (since $c=n!=1\cdot2\cdot\cdots\cdot n$); hence,
$c-1\in\mathbb{N}$. Thus, $n+\left(  c-1\right)  \in\mathbb{N}$ (since
$n\in\mathbb{N}$).

Proposition \ref{prop.sol.choose.a/b.lem1} shows that there exists an
$s\in\mathbb{Z}$ such that $b^{c-1}\equiv sb^{c}\operatorname{mod}c$. Consider
this $s$.

For every $i\in\mathbb{Z}$, we have%
\begin{align}
b^{c-1}\left(  a-bi\right)   &  =\underbrace{b^{c-1}}_{\equiv sb^{c}%
\operatorname{mod}c}a-\underbrace{b^{c-1}b}_{=b^{c}}i\nonumber\\
&  \equiv sb^{c}a-b^{c}i=\underbrace{b^{c}}_{=b^{c-1}b}\left(  sa-i\right)
=b^{c-1}b\left(  sa-i\right)  \operatorname{mod}c.
\label{sol.choose.a/b.congr}%
\end{align}
Thus, in particular, (\ref{sol.choose.a/b.congr}) holds for every
$i\in\left\{  0,1,\ldots,n-1\right\}  $. Hence, Lemma
\ref{lem.sol.choose.a/b.prod-of-congs} (applied to $d=b^{c-1}$, $u_{i}=a-bi$
and $v_{i}=b\left(  sa-i\right)  $) yields%
\begin{equation}
b^{c-1}\prod_{i=0}^{n-1}\left(  a-bi\right)  \equiv b^{c-1}\prod_{i=0}%
^{n-1}\left(  b\left(  sa-i\right)  \right)  \operatorname{mod}c.
\label{sol.choose.a/b.congr.prod}%
\end{equation}

For each $m\in\mathbb{Q}$, we have%
\begin{align}
\dbinom{m}{n}  &  =\dfrac{m\left(  m-1\right)  \cdots\left(  m-n+1\right)
}{n!}=\underbrace{\dfrac{1}{n!}}_{\substack{=\dfrac{1}{c}\\\text{(since
}n!=c\text{)}}}\underbrace{m\left(  m-1\right)  \cdots\left(  m-n+1\right)
}_{=\prod_{i=0}^{n-1}\left(  m-i\right)  }\nonumber\\
&  =\dfrac{1}{c}\prod_{i=0}^{n-1}\left(  m-i\right)  .
\label{sol.choose.a/b.congr.Xn}%
\end{align}
Applying this to $m=sa$, we obtain%
\[
\dbinom{sa}{n}=\dfrac{1}{c}\prod_{i=0}^{n-1}\left(  sa-i\right)  =\dfrac
{\prod_{i=0}^{n-1}\left(  sa-i\right)  }{c}.
\]
Thus,%
\[
\dfrac{\prod_{i=0}^{n-1}\left(  sa-i\right)  }{c}=\dbinom{sa}{n}\in\mathbb{Z}%
\]
(by (\ref{eq.binom.int}) (applied to $m=sa$)). In other words, $c\mid
\prod_{i=0}^{n-1}\left(  sa-i\right)  $. In other words,
\begin{equation}
\prod_{i=0}^{n-1}\left(  sa-i\right)  \equiv0\operatorname{mod}c.
\label{sol.choose.a/b.congr.prod0}%
\end{equation}
Now, (\ref{sol.choose.a/b.congr.prod}) becomes%
\[
b^{c-1}\prod_{i=0}^{n-1}\left(  a-bi\right)  \equiv b^{c-1}\underbrace{\prod
_{i=0}^{n-1}\left(  b\left(  sa-i\right)  \right)  }_{=b^{n}\prod_{i=0}%
^{n-1}\left(  sa-i\right)  }=b^{c-1}b^{n}\underbrace{\prod_{i=0}^{n-1}\left(
sa-i\right)  }_{\substack{\equiv0\operatorname{mod}c\\\text{(by
(\ref{sol.choose.a/b.congr.prod0}))}}}\equiv0\operatorname{mod}c.
\]
In other words,
\begin{equation}
c\mid b^{c-1}\prod_{i=0}^{n-1}\left(  a-bi\right)  .
\label{sol.choose.a/b.divi}%
\end{equation}

But (\ref{sol.choose.a/b.congr.Xn}) (applied to $m=a/b$) yields%
\begin{align*}
\dbinom{a/b}{n}  &  =\dfrac{1}{c}\prod_{i=0}^{n-1}\underbrace{\left(
a/b-i\right)  }_{=\dfrac{1}{b}\left(  a-bi\right)  }=\dfrac{1}{c}%
\underbrace{\prod_{i=0}^{n-1}\left(  \dfrac{1}{b}\left(  a-bi\right)  \right)
}_{=\left(  \dfrac{1}{b}\right)  ^{n}\prod_{i=0}^{n-1}\left(  a-bi\right)
}=\dfrac{1}{c}\underbrace{\left(  \dfrac{1}{b}\right)  ^{n}}_{=\dfrac{1}%
{b^{n}}}\prod_{i=0}^{n-1}\left(  a-bi\right) \\
&  =\dfrac{1}{c}\cdot\dfrac{1}{b^{n}}\prod_{i=0}^{n-1}\left(  a-bi\right)  .
\end{align*}
Multiplying this equality with $b^{n+\left(  c-1\right)  }$, we obtain%
\begin{align*}
b^{n+\left(  c-1\right)  }\dbinom{a/b}{n}  &  =b^{n+\left(  c-1\right)  }%
\cdot\dfrac{1}{c}\cdot\dfrac{1}{b^{n}}\prod_{i=0}^{n-1}\left(  a-bi\right)
=\dfrac{1}{c}\cdot\underbrace{b^{n+\left(  c-1\right)  }\cdot\dfrac{1}{b^{n}}%
}_{=b^{c-1}}\ \ \ \prod_{i=0}^{n-1}\left(  a-bi\right) \\
&  =\dfrac{1}{c}\cdot b^{c-1}\prod_{i=0}^{n-1}\left(  a-bi\right)
=\dfrac{b^{c-1}\prod_{i=0}^{n-1}\left(  a-bi\right)  }{c}\in\mathbb{Z}%
\end{align*}
(by (\ref{sol.choose.a/b.divi})). Hence, there exists some $N\in\mathbb{N}$
such that $b^{N}\dbinom{a/b}{n}\in\mathbb{Z}$ (namely, $N=n+\left(
c-1\right)  $). This solves Exercise \ref{exe.choose.a/b}.
\end{proof}

\begin{remark}
Our above solution of Exercise \ref{exe.choose.a/b} shows that the $N$ in this
exercise can be taken to be $n+\left(  n!-1\right)  $. This is, however, far
from being the best possible value of $N$. A much better value that also works
is $\max\left\{  0,2n-1\right\}  $. Proving this, however, would require a
different idea. The second solution below gives a proof of this better value.
(Alternatively, this better value can be obtained by studying the exponents of
primes appearing in $n!$.)
\end{remark}

\subsubsection{Second solution}

We shall now prepare to give a second solution to Exercise
\ref{exe.choose.a/b}. Our main goal is to prove the following fact:

\begin{theorem}
\label{thm.sol.choose.a/b.sol2.main}Let $a$ and $b$ be two integers such that
$b\neq0$. Let $n$ be a positive integer. Then, $b^{2n-1}\dbinom{a/b}{n}%
\in\mathbb{Z}$.
\end{theorem}

Clearly, Exercise \ref{exe.choose.a/b} immediately follows from Theorem
\ref{thm.sol.choose.a/b.sol2.main} in the case when $n\neq0$. (In the case
when $n=0$, it holds for obvious reasons.)

Before we start proving Theorem \ref{thm.sol.choose.a/b.sol2.main}, let us
show the following lemma:

\begin{lemma}
\label{lem.sol.choose.a/b.sol2.lem}Let $b$ and $n$ be positive integers.
Assume that every $k\in\left\{  1,2,\ldots,n-1\right\}  $ and $c\in\mathbb{Z}$
satisfy%
\begin{equation}
b^{2k-1}\dbinom{c/b}{k}\in\mathbb{Z}.
\label{eq.lem.sol.choose.a/b.sol2.lem.ass}%
\end{equation}
Then:

\textbf{(a)} Every $u\in\mathbb{Z}$ and $v\in\mathbb{Z}$ satisfy%
\[
b^{2n-2}\left(  \dbinom{\left(  u+v\right)  /b}{n}-\dbinom{u/b}{n}%
-\dbinom{v/b}{n}\right)  \in\mathbb{Z}.
\]

\textbf{(b)} Every $a\in\mathbb{Z}$ and $h\in\mathbb{N}$ satisfy%
\[
b^{2n-2}\left(  \dbinom{ha/b}{n}-h\dbinom{a/b}{n}\right)  \in\mathbb{Z}.
\]

\textbf{(c)} Every $a\in\mathbb{Z}$ satisfies%
\[
b^{2n-1}\dbinom{a/b}{n}\in\mathbb{Z}.
\]

\end{lemma}

\begin{vershort}
\begin{proof}
[Proof of Lemma \ref{lem.sol.choose.a/b.sol2.lem}.]\textbf{(a)} Let
$u\in\mathbb{Z}$ and $v\in\mathbb{Z}$. But Proposition \ref{prop.binom.00}
\textbf{(a)} (applied to $u/b$ instead of $m$) yields $\dbinom{u/b}{0}=1$.
Similarly, $\dbinom{v/b}{0}=1$.

Theorem \ref{thm.vandermonde.rat} (applied to $u/b$ and $v/b$ instead of $x$
and $y$) yields
\begin{align*}
\dbinom{u/b+v/b}{n}  &  =\sum_{k=0}^{n}\dbinom{u/b}{k}\dbinom{v/b}{n-k}\\
&  =\underbrace{\dbinom{u/b}{0}}_{=1}\underbrace{\dbinom{v/b}{n-0}}%
_{=\dbinom{v/b}{n}}+\sum_{k=1}^{n-1}\dbinom{u/b}{k}\dbinom{v/b}{n-k}%
+\dbinom{u/b}{n}\underbrace{\dbinom{v/b}{n-n}}_{=\dbinom{v/b}{0}=1}\\
&  \ \ \ \ \ \ \ \ \ \ \left(
\begin{array}
[c]{c}%
\text{here, we have split off the addends for }k=0\text{ and}\\
\text{for }k=n\text{ from the sum (since }0\text{ and }n\text{ are two
distinct}\\
\text{elements of }\left\{  0,1,\ldots,n\right\}  \text{)}%
\end{array}
\right) \\
&  =\dbinom{v/b}{n}+\sum_{k=1}^{n-1}\dbinom{u/b}{k}\dbinom{v/b}{n-k}%
+\dbinom{u/b}{n}.
\end{align*}
Subtracting $\dbinom{u/b}{n}+\dbinom{v/b}{n}$ from this equality, we obtain%
\begin{equation}
\dbinom{u/b+v/b}{n}-\dbinom{u/b}{n}-\dbinom{v/b}{n}=\sum_{k=1}^{n-1}%
\dbinom{u/b}{k}\dbinom{v/b}{n-k}.
\label{pf.lem.sol.choose.a/b.sol2.lem.short.a.1}%
\end{equation}

But for every $k\in\left\{  1,2,\ldots,n-1\right\}  $, the number
$b^{2n-2}\dbinom{u/b}{k}\dbinom{v/b}{n-k}$ is an
integer\footnote{\textit{Proof:} Let $k\in\left\{  1,2,\ldots,n-1\right\}  $.
\par
From (\ref{eq.lem.sol.choose.a/b.sol2.lem.ass}) (applied to $c=u$), we obtain
$b^{2k-1}\dbinom{u/b}{k}\in\mathbb{Z}$. In other words, the number
$b^{2k-1}\dbinom{u/b}{k}$ is an integer.
\par
But we also have $n-k\in\left\{  1,2,\ldots,n-1\right\}  $ (since
$k\in\left\{  1,2,\ldots,n-1\right\}  $). Thus, from
(\ref{eq.lem.sol.choose.a/b.sol2.lem.ass}) (applied to $n-k$ and $v$ instead
of $k$ and $c$), we obtain $b^{2\left(  n-k\right)  -1}\dbinom{v/b}{n-k}%
\in\mathbb{Z}$. In other words, the number $b^{2\left(  n-k\right)  -1}%
\dbinom{v/b}{n-k}$ is an integer.
\par
Now, the numbers $b^{2k-1}\dbinom{u/b}{k}$ and $b^{2\left(  n-k\right)
-1}\dbinom{v/b}{n-k}$ are integers. Hence, their product is an integer as
well. In other words, the number $b^{2k-1}\dbinom{u/b}{k}\cdot b^{2\left(
n-k\right)  -1}\dbinom{v/b}{n-k}$ is an integer. Since%
\begin{align*}
&  b^{2k-1}\dbinom{u/b}{k}\cdot b^{2\left(  n-k\right)  -1}\dbinom{v/b}{n-k}\\
&  =\underbrace{b^{2k-1}b^{2\left(  n-k\right)  -1}}_{\substack{=b^{\left(
2k-1\right)  +\left(  2\left(  n-k\right)  -1\right)  }=b^{2n-2}\\\text{(since
}\left(  2k-1\right)  +\left(  2\left(  n-k\right)  -1\right)  =2n-2\text{)}%
}}\dbinom{u/b}{k}\dbinom{v/b}{n-k}=b^{2n-2}\dbinom{u/b}{k}\dbinom{v/b}{n-k},
\end{align*}
this rewrites as follows: The number $b^{2n-2}\dbinom{u/b}{k}\dbinom{v/b}%
{n-k}$ is an integer. Qed.}. Hence, $\sum_{k=1}^{n-1}b^{2n-2}\dbinom{u/b}%
{k}\dbinom{v/b}{n-k}$ is a sum of $n-1$ integers, and thus itself is an
integer. In other words,%
\begin{equation}
\sum_{k=1}^{n-1}b^{2n-2}\dbinom{u/b}{k}\dbinom{v/b}{n-k}\in\mathbb{Z}.
\label{pf.lem.sol.choose.a/b.sol2.lem.short.a.5}%
\end{equation}

Now,%
\begin{align*}
&  b^{2n-2}\left(  \underbrace{\dbinom{\left(  u+v\right)  /b}{n}%
}_{\substack{=\dbinom{u/b+v/b}{n}\\\text{(since }\left(  u+v\right)
/b=u/b+v/b\text{)}}}-\dbinom{u/b}{n}-\dbinom{v/b}{n}\right) \\
&  =b^{2n-2}\underbrace{\left(  \dbinom{u/b+v/b}{n}-\dbinom{u/b}{n}%
-\dbinom{v/b}{n}\right)  }_{\substack{=\sum_{k=1}^{n-1}\dbinom{u/b}{k}%
\dbinom{v/b}{n-k}\\\text{(by (\ref{pf.lem.sol.choose.a/b.sol2.lem.short.a.1}%
))}}}\\
&  =b^{2n-2}\sum_{k=1}^{n-1}\dbinom{u/b}{k}\dbinom{v/b}{n-k}=\sum_{k=1}%
^{n-1}b^{2n-2}\dbinom{u/b}{k}\dbinom{v/b}{n-k}\in\mathbb{Z}%
\end{align*}
(by (\ref{pf.lem.sol.choose.a/b.sol2.lem.short.a.5})). This proves Lemma
\ref{lem.sol.choose.a/b.sol2.lem} \textbf{(a)}.

\textbf{(b)} Let $a\in\mathbb{Z}$ and $h\in\mathbb{N}$. For every
$k\in\left\{  0,1,\ldots,h-1\right\}  $, the number \newline$b^{2n-2}\left(
\dbinom{\left(  k+1\right)  a/b}{n}-\dbinom{ka/b}{n}-\dbinom{a/b}{n}\right)  $
is an integer\footnote{\textit{Proof:} Let $k\in\left\{  0,1,\ldots
,h-1\right\}  $. Then, Lemma \ref{lem.sol.choose.a/b.sol2.lem} \textbf{(a)}
(applied to $u=ka$ and $v=a$) yields%
\[
b^{2n-2}\left(  \dbinom{\left(  ka+a\right)  /b}{n}-\dbinom{ka/b}{n}%
-\dbinom{a/b}{n}\right)  \in\mathbb{Z}.
\]
Since $ka+a=\left(  k+1\right)  a$, this rewrites as follows:%
\[
b^{2n-2}\left(  \dbinom{\left(  k+1\right)  a/b}{n}-\dbinom{ka/b}{n}%
-\dbinom{a/b}{n}\right)  \in\mathbb{Z}.
\]
Qed.}. Hence,
\[
\sum_{k=0}^{h-1}b^{2n-2}\left(  \dbinom{\left(  k+1\right)  a/b}{n}%
-\dbinom{ka/b}{n}-\dbinom{a/b}{n}\right)
\]
is a sum of $h$ integers, and thus itself is an integer. In other words,%
\begin{equation}
\sum_{k=0}^{h-1}b^{2n-2}\left(  \dbinom{\left(  k+1\right)  a/b}{n}%
-\dbinom{ka/b}{n}-\dbinom{a/b}{n}\right)  \in\mathbb{Z}.
\label{pf.lem.sol.choose.a/b.sol2.lem.short.b.3}%
\end{equation}

But $1-1=0\leq h$. Hence, (\ref{eq.sum.telescope}) (applied to $\mathbb{A}%
=\mathbb{Q}$, $u=1$, $v=h$ and $a_{s}=\dbinom{sa/b}{n}$) yields%
\begin{align*}
\sum_{s=1}^{h}\left(  \dbinom{sa/b}{n}-\dbinom{\left(  s-1\right)  a/b}%
{n}\right)   &  =\dbinom{ha/b}{n}-\underbrace{\dbinom{\left(  1-1\right)
a/b}{n}}_{\substack{=\dbinom{0}{n}\\\text{(since }\left(  1-1\right)
a/b=0\text{)}}}\\
&  =\dbinom{ha/b}{n}-\underbrace{\dbinom{0}{n}}_{\substack{=0\\\text{(by
Proposition \ref{prop.binom.0},}\\\text{applied to }m=0\text{)}}%
}=\dbinom{ha/b}{n}.
\end{align*}
Hence,%
\begin{align}
\dbinom{ha/b}{n}  &  =\sum_{s=1}^{h}\left(  \dbinom{sa/b}{n}-\dbinom{\left(
s-1\right)  a/b}{n}\right) \nonumber\\
&  =\sum_{k=0}^{h-1}\left(  \dbinom{\left(  k+1\right)  a/b}{n}-\dbinom
{ka/b}{n}\right)  \label{pf.lem.sol.choose.a/b.sol2.lem.short.b.6}%
\end{align}
(here, we have substituted $k+1$ for $s$ in the sum). Now,%
\begin{align}
&  \sum_{k=0}^{h-1}\left(  \dbinom{\left(  k+1\right)  a/b}{n}-\dbinom
{ka/b}{n}-\dbinom{a/b}{n}\right) \nonumber\\
&  =\underbrace{\sum_{k=0}^{h-1}\left(  \dbinom{\left(  k+1\right)  a/b}%
{n}-\dbinom{ka/b}{n}\right)  }_{\substack{=\dbinom{ha/b}{n}\\\text{(by
(\ref{pf.lem.sol.choose.a/b.sol2.lem.short.b.6}))}}}-\underbrace{\sum
_{k=0}^{h-1}\dbinom{a/b}{n}}_{=h\dbinom{a/b}{n}}\nonumber\\
&  =\dbinom{ha/b}{n}-h\dbinom{a/b}{n}.
\label{pf.lem.sol.choose.a/b.sol2.lem.short.b.7}%
\end{align}

Now,
\begin{align*}
&  \sum_{k=0}^{h-1}b^{2n-2}\left(  \dbinom{\left(  k+1\right)  a/b}{n}%
-\dbinom{ka/b}{n}-\dbinom{a/b}{n}\right) \\
&  =b^{2n-2}\underbrace{\sum_{k=0}^{h-1}\left(  \dbinom{\left(  k+1\right)
a/b}{n}-\dbinom{ka/b}{n}-\dbinom{a/b}{n}\right)  }_{\substack{=\dbinom
{ha/b}{n}-h\dbinom{a/b}{n}\\\text{(by
(\ref{pf.lem.sol.choose.a/b.sol2.lem.short.b.7}))}}}\\
&  =b^{2n-2}\left(  \dbinom{ha/b}{n}-h\dbinom{a/b}{n}\right)  .
\end{align*}
Thus,%
\begin{align*}
b^{2n-2}\left(  \dbinom{ha/b}{n}-h\dbinom{a/b}{n}\right)   &  =\sum
_{k=0}^{h-1}b^{2n-2}\left(  \dbinom{\left(  k+1\right)  a/b}{n}-\dbinom
{ka/b}{n}-\dbinom{a/b}{n}\right) \\
&  \in\mathbb{Z}\ \ \ \ \ \ \ \ \ \ \left(  \text{by
(\ref{pf.lem.sol.choose.a/b.sol2.lem.short.b.3})}\right)  .
\end{align*}
This proves Lemma \ref{lem.sol.choose.a/b.sol2.lem} \textbf{(b)}.

\textbf{(c)} Let $a\in\mathbb{Z}$. Proposition \ref{prop.binom.int} (applied
to $m=a$) yields $\dbinom{a}{n}\in\mathbb{Z}$. In other words, $\dbinom{a}{n}$
is an integer.

Also, $2n-2\in\mathbb{N}$ (since $n$ is a positive integer). Thus, $b^{2n-2}$
is an integer (since $b$ is an integer). Now, the numbers $b^{2n-2}$ and
$\dbinom{a}{n}$ are both integers. Hence, their product must also be an
integer. In other words, $b^{2n-2}\dbinom{a}{n}$ is an integer.

But Lemma \ref{lem.sol.choose.a/b.sol2.lem} \textbf{(b)} (applied to $h=b$)
yields%
\[
b^{2n-2}\left(  \dbinom{ba/b}{n}-b\dbinom{a/b}{n}\right)  \in\mathbb{Z}.
\]
In other words, $b^{2n-2}\left(  \dbinom{ba/b}{n}-b\dbinom{a/b}{n}\right)  $
is an integer. Denote this integer by $z$. Thus,%
\begin{align*}
z  &  =b^{2n-2}\left(  \dbinom{ba/b}{n}-b\dbinom{a/b}{n}\right)
=b^{2n-2}\underbrace{\dbinom{ba/b}{n}}_{=\dbinom{a}{n}}-\underbrace{b^{2n-2}%
b}_{=b^{\left(  2n-2\right)  +1}=b^{2n-1}}\dbinom{a/b}{n}\\
&  =b^{2n-2}\dbinom{a}{n}-b^{2n-1}\dbinom{a/b}{n}.
\end{align*}
Solving this equality for $b^{2n-1}\dbinom{a/b}{n}$, we obtain
\begin{equation}
b^{2n-1}\dbinom{a/b}{n}=b^{2n-2}\dbinom{a}{n}-z.
\label{pf.lem.sol.choose.a/b.sol2.lem.short.c.3}%
\end{equation}

But the numbers $b^{2n-2}\dbinom{a}{n}$ and $z$ are integers. Hence, their
difference is also an integer. In other words, $b^{2n-2}\dbinom{a}{n}-z$ is an
integer. In other words, $b^{2n-2}\dbinom{a}{n}-z\in\mathbb{Z}$. Hence,
(\ref{pf.lem.sol.choose.a/b.sol2.lem.short.c.3}) becomes $b^{2n-1}\dbinom
{a/b}{n}=b^{2n-2}\dbinom{a}{n}-z\in\mathbb{Z}$. This proves Lemma
\ref{lem.sol.choose.a/b.sol2.lem} \textbf{(c)}.
\end{proof}
\end{vershort}

\begin{verlong}
\begin{proof}
[Proof of Lemma \ref{lem.sol.choose.a/b.sol2.lem}.]The integer $n$ is
positive. Hence, $n-1\in\mathbb{N}$, so that $0\in\left\{  0,1,\ldots
,n-1\right\}  $. Also, $n\in\mathbb{N}$ (since $n$ is a positive integer), so
that $n\in\left\{  0,1,\ldots,n\right\}  $.

\textbf{(a)} Let $u\in\mathbb{Z}$ and $v\in\mathbb{Z}$. But Proposition
\ref{prop.binom.00} \textbf{(a)} (applied to $u/b$ instead of $m$) yields
$\dbinom{u/b}{0}=1$. Also, Proposition \ref{prop.binom.00} \textbf{(a)}
(applied to $v/b$ instead of $m$) yields $\dbinom{v/b}{0}=1$.

Theorem \ref{thm.vandermonde.rat} (applied to $u/b$ and $v/b$ instead of $x$
and $y$) yields
\begin{align*}
\dbinom{u/b+v/b}{n}  &  =\sum_{k=0}^{n}\dbinom{u/b}{k}\dbinom{v/b}{n-k}\\
&  =\sum_{k=0}^{n-1}\dbinom{u/b}{k}\dbinom{v/b}{n-k}+\dbinom{u/b}%
{n}\underbrace{\dbinom{v/b}{n-n}}_{=\dbinom{v/b}{0}=1}\\
&  \ \ \ \ \ \ \ \ \ \ \left(
\begin{array}
[c]{c}%
\text{here, we have split off the addend for }k=n\text{ from}\\
\text{the sum (since }n\in\left\{  0,1,\ldots,n\right\}  \text{)}%
\end{array}
\right) \\
&  =\sum_{k=0}^{n-1}\dbinom{u/b}{k}\dbinom{v/b}{n-k}+\dbinom{u/b}{n}.
\end{align*}
Subtracting $\dbinom{u/b}{n}$ from this equality, we obtain%
\begin{align*}
&  \dbinom{u/b+v/b}{n}-\dbinom{u/b}{n}\\
&  =\sum_{k=0}^{n-1}\dbinom{u/b}{k}\dbinom{v/b}{n-k}\\
&  =\underbrace{\dbinom{u/b}{0}}_{=1}\underbrace{\dbinom{v/b}{n-0}}%
_{=\dbinom{v/b}{n}}+\sum_{k=1}^{n-1}\dbinom{u/b}{k}\dbinom{v/b}{n-k}\\
&  =\dbinom{v/b}{n}+\sum_{k=1}^{n-1}\dbinom{u/b}{k}\dbinom{v/b}{n-k}.
\end{align*}
Subtracting $\dbinom{v/b}{n}$ from this equality, we obtain%
\begin{equation}
\dbinom{u/b+v/b}{n}-\dbinom{u/b}{n}-\dbinom{v/b}{n}=\sum_{k=1}^{n-1}%
\dbinom{u/b}{k}\dbinom{v/b}{n-k}. \label{pf.lem.sol.choose.a/b.sol2.lem.a.1}%
\end{equation}

But for every $k\in\left\{  1,2,\ldots,n-1\right\}  $, we have%
\begin{equation}
b^{2n-2}\dbinom{u/b}{k}\dbinom{v/b}{n-k}\in\mathbb{Z}
\label{pf.lem.sol.choose.a/b.sol2.lem.a.3}%
\end{equation}
\footnote{\textit{Proof of (\ref{pf.lem.sol.choose.a/b.sol2.lem.a.3}):} Let
$k\in\left\{  1,2,\ldots,n-1\right\}  $.
\par
From (\ref{eq.lem.sol.choose.a/b.sol2.lem.ass}) (applied to $c=u$), we obtain
$b^{2k-1}\dbinom{u/b}{k}\in\mathbb{Z}$. In other words, the number
$b^{2k-1}\dbinom{u/b}{k}$ is an integer.
\par
But we also have $n-k\in\left\{  1,2,\ldots,n-1\right\}  $ (since
$k\in\left\{  1,2,\ldots,n-1\right\}  $). Thus, from
(\ref{eq.lem.sol.choose.a/b.sol2.lem.ass}) (applied to $n-k$ and $v$ instead
of $k$ and $c$), we obtain $b^{2\left(  n-k\right)  -1}\dbinom{v/b}{n-k}%
\in\mathbb{Z}$. In other words, the number $b^{2\left(  n-k\right)  -1}%
\dbinom{v/b}{n-k}$ is an integer.
\par
Now, the numbers $b^{2k-1}\dbinom{u/b}{k}$ and $b^{2\left(  n-k\right)
-1}\dbinom{v/b}{n-k}$ are integers. Hence, their product is an integer as
well. In other words, $b^{2k-1}\dbinom{u/b}{k}\cdot b^{2\left(  n-k\right)
-1}\dbinom{v/b}{n-k}$ is an integer. In other words, $b^{2k-1}\dbinom{u/b}%
{k}\cdot b^{2\left(  n-k\right)  -1}\dbinom{v/b}{n-k}\in\mathbb{Z}$. Since%
\begin{align*}
&  b^{2k-1}\dbinom{u/b}{k}\cdot b^{2\left(  n-k\right)  -1}\dbinom{v/b}{n-k}\\
&  =\underbrace{b^{2k-1}b^{2\left(  n-k\right)  -1}}_{\substack{=b^{\left(
2k-1\right)  +\left(  2\left(  n-k\right)  -1\right)  }=b^{2n-2}\\\text{(since
}\left(  2k-1\right)  +\left(  2\left(  n-k\right)  -1\right)  =2n-2\text{)}%
}}\dbinom{u/b}{k}\dbinom{v/b}{n-k}=b^{2n-2}\dbinom{u/b}{k}\dbinom{v/b}{n-k},
\end{align*}
this rewrites as $b^{2n-2}\dbinom{u/b}{k}\dbinom{v/b}{n-k}\in\mathbb{Z}$.
Thus, (\ref{pf.lem.sol.choose.a/b.sol2.lem.a.3}) is proven.}. In other words,
for every $k\in\left\{  1,2,\ldots,n-1\right\}  $, the number $b^{2n-2}%
\dbinom{u/b}{k}\dbinom{v/b}{n-k}$ is an integer. Hence, $\sum_{k=1}%
^{n-1}b^{2n-2}\dbinom{u/b}{k}\dbinom{v/b}{n-k}$ is a sum of $n-1$ integers,
and thus itself is an integer. In other words,%
\begin{equation}
\sum_{k=1}^{n-1}b^{2n-2}\dbinom{u/b}{k}\dbinom{v/b}{n-k}\in\mathbb{Z}.
\label{pf.lem.sol.choose.a/b.sol2.lem.a.5}%
\end{equation}

Now,%
\begin{align*}
&  b^{2n-2}\left(  \underbrace{\dbinom{\left(  u+v\right)  /b}{n}%
}_{\substack{=\dbinom{u/b+v/b}{n}\\\text{(since }\left(  u+v\right)
/b=u/b+v/b\text{)}}}-\dbinom{u/b}{n}-\dbinom{v/b}{n}\right) \\
&  =b^{2n-2}\underbrace{\left(  \dbinom{u/b+v/b}{n}-\dbinom{u/b}{n}%
-\dbinom{v/b}{n}\right)  }_{\substack{=\sum_{k=1}^{n-1}\dbinom{u/b}{k}%
\dbinom{v/b}{n-k}\\\text{(by (\ref{pf.lem.sol.choose.a/b.sol2.lem.a.1}))}}}\\
&  =b^{2n-2}\sum_{k=1}^{n-1}\dbinom{u/b}{k}\dbinom{v/b}{n-k}=\sum_{k=1}%
^{n-1}b^{2n-2}\dbinom{u/b}{k}\dbinom{v/b}{n-k}\in\mathbb{Z}%
\end{align*}
(by (\ref{pf.lem.sol.choose.a/b.sol2.lem.a.5})). This proves Lemma
\ref{lem.sol.choose.a/b.sol2.lem} \textbf{(a)}.

\textbf{(b)} Let $a\in\mathbb{Z}$ and $h\in\mathbb{N}$. For every
$k\in\left\{  0,1,\ldots,h-1\right\}  $, we have%
\begin{equation}
b^{2n-2}\left(  \dbinom{\left(  k+1\right)  a/b}{n}-\dbinom{ka/b}{n}%
-\dbinom{a/b}{n}\right)  \in\mathbb{Z}
\label{pf.lem.sol.choose.a/b.sol2.lem.b.1}%
\end{equation}
\footnote{\textit{Proof of (\ref{pf.lem.sol.choose.a/b.sol2.lem.b.1}):} Let
$k\in\left\{  0,1,\ldots,h-1\right\}  $. Then, Lemma
\ref{lem.sol.choose.a/b.sol2.lem} \textbf{(a)} (applied to $u=ka$ and $v=a$)
yields%
\[
b^{2n-2}\left(  \dbinom{\left(  ka+a\right)  /b}{n}-\dbinom{ka/b}{n}%
-\dbinom{a/b}{n}\right)  \in\mathbb{Z}.
\]
Since $ka+a=\left(  k+1\right)  a$, this rewrites as follows:%
\[
b^{2n-2}\left(  \dbinom{\left(  k+1\right)  a/b}{n}-\dbinom{ka/b}{n}%
-\dbinom{a/b}{n}\right)  \in\mathbb{Z}.
\]
Thus, (\ref{pf.lem.sol.choose.a/b.sol2.lem.b.1}) is proven.}. In other words,
for every $k\in\left\{  0,1,\ldots,h-1\right\}  $, the number \newline%
$b^{2n-2}\left(  \dbinom{\left(  k+1\right)  a/b}{n}-\dbinom{ka/b}{n}%
-\dbinom{a/b}{n}\right)  $ is an integer. Hence,
\[
\sum_{k=0}^{h-1}b^{2n-2}\left(  \dbinom{\left(  k+1\right)  a/b}{n}%
-\dbinom{ka/b}{n}-\dbinom{a/b}{n}\right)
\]
is a sum of $h$ integers, and thus itself is an integer. In other words,%
\begin{equation}
\sum_{k=0}^{h-1}b^{2n-2}\left(  \dbinom{\left(  k+1\right)  a/b}{n}%
-\dbinom{ka/b}{n}-\dbinom{a/b}{n}\right)  \in\mathbb{Z}.
\label{pf.lem.sol.choose.a/b.sol2.lem.b.3}%
\end{equation}

But $h\in\mathbb{N}$ and thus $h\in\left\{  0,1,\ldots,h\right\}  $. Hence,
\begin{align*}
\sum_{k=0}^{h}\dbinom{ka/b}{n}  &  =\sum_{k=0}^{h-1}\dbinom{ka/b}{n}%
+\dbinom{ha/b}{n}\\
&  \ \ \ \ \ \ \ \ \ \ \left(
\begin{array}
[c]{c}%
\text{here, we have split off the addend for }k=h\\
\text{from the sum (since }h\in\left\{  0,1,\ldots,h\right\}  \text{)}%
\end{array}
\right)  .
\end{align*}
Subtracting $\dbinom{ha/b}{n}$ from this equality, we obtain%
\begin{equation}
\sum_{k=0}^{h}\dbinom{ka/b}{n}-\dbinom{ha/b}{n}=\sum_{k=0}^{h-1}\dbinom
{ka/b}{n}. \label{pf.lem.sol.choose.a/b.sol2.lem.b.5a}%
\end{equation}

On the other hand, $n$ is a positive integer. Hence, $0<n$. Thus, Proposition
\ref{prop.binom.0} (applied to $m=0$) yields $\dbinom{0}{n}=0$.

But $h\in\mathbb{N}$ and thus $0\in\left\{  0,1,\ldots,h\right\}  $. Hence,
\begin{align}
\sum_{k=0}^{h}\dbinom{ka/b}{n}  &  =\underbrace{\dbinom{0a/b}{n}%
}_{\substack{=\dbinom{0}{n}\\\text{(since }0a/b=0\text{)}}}+\sum_{k=1}%
^{h}\dbinom{ka/b}{n}\nonumber\\
&  \ \ \ \ \ \ \ \ \ \ \left(
\begin{array}
[c]{c}%
\text{here, we have split off the addend for }k=0\\
\text{from the sum (since }0\in\left\{  0,1,\ldots,h\right\}  \text{)}%
\end{array}
\right) \nonumber\\
&  =\underbrace{\dbinom{0}{n}}_{=0}+\sum_{k=1}^{h}\dbinom{ka/b}{n}=\sum
_{k=1}^{h}\dbinom{ka/b}{n}\nonumber\\
&  =\sum_{k=0}^{h-1}\dbinom{\left(  k+1\right)  a/b}{n}
\label{pf.lem.sol.choose.a/b.sol2.lem.b.5b}%
\end{align}
(here, we have substituted $k+1$ for $k$ in the sum).

But%
\begin{align}
&  \sum_{k=0}^{h-1}\left(  \dbinom{\left(  k+1\right)  a/b}{n}-\dbinom
{ka/b}{n}-\dbinom{a/b}{n}\right) \nonumber\\
&  =\underbrace{\sum_{k=0}^{h-1}\dbinom{\left(  k+1\right)  a/b}{n}%
}_{\substack{=\sum_{k=0}^{h}\dbinom{ka/b}{n}\\\text{(by
(\ref{pf.lem.sol.choose.a/b.sol2.lem.b.5b}))}}}-\underbrace{\sum_{k=0}%
^{h-1}\dbinom{ka/b}{n}}_{\substack{=\sum_{k=0}^{h}\dbinom{ka/b}{n}%
-\dbinom{ha/b}{n}\\\text{(by (\ref{pf.lem.sol.choose.a/b.sol2.lem.b.5a}))}%
}}-\underbrace{\sum_{k=0}^{h-1}\dbinom{a/b}{n}}_{=h\dbinom{a/b}{n}}\nonumber\\
&  =\sum_{k=0}^{h}\dbinom{ka/b}{n}-\left(  \sum_{k=0}^{h}\dbinom{ka/b}%
{n}-\dbinom{ha/b}{n}\right)  -h\dbinom{a/b}{n}\nonumber\\
&  =\dbinom{ha/b}{n}-h\dbinom{a/b}{n}.
\label{pf.lem.sol.choose.a/b.sol2.lem.b.7}%
\end{align}

Now,
\begin{align*}
&  \sum_{k=0}^{h-1}b^{2n-2}\left(  \dbinom{\left(  k+1\right)  a/b}{n}%
-\dbinom{ka/b}{n}-\dbinom{a/b}{n}\right) \\
&  =b^{2n-2}\underbrace{\sum_{k=0}^{h-1}\left(  \dbinom{\left(  k+1\right)
a/b}{n}-\dbinom{ka/b}{n}-\dbinom{a/b}{n}\right)  }_{\substack{=\dbinom
{ha/b}{n}-h\dbinom{a/b}{n}\\\text{(by
(\ref{pf.lem.sol.choose.a/b.sol2.lem.b.7}))}}}\\
&  =b^{2n-2}\left(  \dbinom{ha/b}{n}-h\dbinom{a/b}{n}\right)  .
\end{align*}
Thus,%
\begin{align*}
&  b^{2n-2}\left(  \dbinom{ha/b}{n}-h\dbinom{a/b}{n}\right) \\
&  =\sum_{k=0}^{h-1}b^{2n-2}\left(  \dbinom{\left(  k+1\right)  a/b}%
{n}-\dbinom{ka/b}{n}-\dbinom{a/b}{n}\right) \\
&  \in\mathbb{Z}\ \ \ \ \ \ \ \ \ \ \left(  \text{by
(\ref{pf.lem.sol.choose.a/b.sol2.lem.b.3})}\right)  .
\end{align*}
This proves Lemma \ref{lem.sol.choose.a/b.sol2.lem} \textbf{(b)}.

\textbf{(c)} Let $a\in\mathbb{Z}$. Proposition \ref{prop.binom.int} (applied
to $m=a$) yields $\dbinom{a}{n}\in\mathbb{Z}$. In other words, $\dbinom{a}{n}$
is an integer.

Also, $n$ is a positive integer; thus, $n\geq1$, so that $n-1\in\mathbb{N}$
and thus $2\left(  n-1\right)  \in\mathbb{N}$. Hence, $2n-2=2\left(
n-1\right)  \in\mathbb{N}$. Thus, $b^{2n-2}$ is an integer (since $b$ is an
integer). Now, the numbers $b^{2n-2}$ and $\dbinom{a}{n}$ are both integers.
Hence, their product must also be an integer. In other words, $b^{2n-2}%
\dbinom{a}{n}$ is an integer.

But Lemma \ref{lem.sol.choose.a/b.sol2.lem} \textbf{(b)} (applied to $h=b$)
yields%
\[
b^{2n-2}\left(  \dbinom{ba/b}{n}-b\dbinom{a/b}{n}\right)  \in\mathbb{Z}.
\]
In other words, $b^{2n-2}\left(  \dbinom{ba/b}{n}-b\dbinom{a/b}{n}\right)  $
is an integer. Denote this integer by $z$. Thus,%
\begin{align*}
z  &  =b^{2n-2}\left(  \dbinom{ba/b}{n}-b\dbinom{a/b}{n}\right) \\
&  =b^{2n-2}\underbrace{\dbinom{ba/b}{n}}_{\substack{=\dbinom{a}%
{n}\\\text{(since }ba/b=a\text{)}}}-\underbrace{b^{2n-2}b}%
_{\substack{=b^{\left(  2n-2\right)  +1}=b^{2n-1}\\\text{(since }\left(
2n-2\right)  +1=2n-1\text{)}}}\dbinom{a/b}{n}\\
&  =b^{2n-2}\dbinom{a}{n}-b^{2n-1}\dbinom{a/b}{n}.
\end{align*}
Adding $b^{2n-1}\dbinom{a/b}{n}$ to both sides of this equality, we obtain
$b^{2n-1}\dbinom{a/b}{n}+z=b^{2n-2}\dbinom{a}{n}$. Subtracting $z$ from this
equality, we obtain
\begin{equation}
b^{2n-1}\dbinom{a/b}{n}=b^{2n-2}\dbinom{a}{n}-z.
\label{pf.lem.sol.choose.a/b.sol2.lem.c.3}%
\end{equation}

But the numbers $b^{2n-2}\dbinom{a}{n}$ and $z$ are integers. Hence, their
difference is also an integer. In other words, $b^{2n-2}\dbinom{a}{n}-z$ is an
integer. In other words, $b^{2n-2}\dbinom{a}{n}-z\in\mathbb{Z}$. Hence,
(\ref{pf.lem.sol.choose.a/b.sol2.lem.c.3}) becomes $b^{2n-1}\dbinom{a/b}%
{n}=b^{2n-2}\dbinom{a}{n}-z\in\mathbb{Z}$. This proves Lemma
\ref{lem.sol.choose.a/b.sol2.lem} \textbf{(c)}.
\end{proof}
\end{verlong}

Our next lemma is essentially Theorem \ref{thm.sol.choose.a/b.sol2.main},
restricted to the case when $b$ is positive:

\begin{lemma}
\label{lem.sol.choose.a/b.sol2.pos}Let $a$ and $b$ be two integers such that
$b>0$. Let $n$ be a positive integer. Then, $b^{2n-1}\dbinom{a/b}{n}%
\in\mathbb{Z}$.
\end{lemma}

\begin{proof}
[Proof of Lemma \ref{lem.sol.choose.a/b.sol2.pos}.]We shall prove Lemma
\ref{lem.sol.choose.a/b.sol2.pos} by strong induction on $n$:

\textit{Induction step:} Let $N$ be a positive integer. Assume that Lemma
\ref{lem.sol.choose.a/b.sol2.pos} holds in the case when $n<N$. We must show
that Lemma \ref{lem.sol.choose.a/b.sol2.pos} holds in the case when $n=N$.

We have assumed that Lemma \ref{lem.sol.choose.a/b.sol2.pos} holds in the case
when $n<N$. In other words, the following statement holds:

\begin{statement}
\textit{Statement 1:} Let $a$ and $b$ be two integers such that $b>0$. Let $n$
be a positive integer such that $n<N$. Then, $b^{2n-1}\dbinom{a/b}{n}%
\in\mathbb{Z}$.
\end{statement}

Now, let us prove the following statement:

\begin{statement}
\textit{Statement 2:} Let $a$ and $b$ be two integers such that $b>0$. Then,
$b^{2N-1}\dbinom{a/b}{N}\in\mathbb{Z}$.
\end{statement}

[\textit{Proof of Statement 2:} Every $k\in\left\{  1,2,\ldots,N-1\right\}  $
and $c\in\mathbb{Z}$ satisfy%
\[
b^{2k-1}\dbinom{c/b}{k}\in\mathbb{Z}%
\]
\footnote{\textit{Proof.} Let $k\in\left\{  1,2,\ldots,N-1\right\}  $ and
$c\in\mathbb{Z}$. From $k\in\left\{  1,2,\ldots,N-1\right\}  $, we obtain
$1\leq k\leq N-1$. Now, $k$ is a positive integer (since $1\leq k$) and
satisfies $k<N$ (since $k\leq N-1<N$). Hence, Statement 1 (applied to $k$ and
$c$ instead of $n$ and $a$) yields $b^{2k-1}\dbinom{c/b}{k}\in\mathbb{Z}$.
Qed.}. Hence, Lemma \ref{lem.sol.choose.a/b.sol2.lem} \textbf{(c)} (applied to
$n=N$) yields $b^{2N-1}\dbinom{a/b}{N}\in\mathbb{Z}$. Thus, Statement 2 is proven.]

So we have proven Statement 2. In other words, we have proven that Lemma
\ref{lem.sol.choose.a/b.sol2.pos} holds in the case when $n=N$. This completes
the induction step. Thus, the inductive proof of Lemma
\ref{lem.sol.choose.a/b.sol2.pos} is complete.
\end{proof}

\begin{vershort}
\begin{proof}
[Proof of Theorem \ref{thm.sol.choose.a/b.sol2.main}.]We must prove that
$b^{2n-1}\dbinom{a/b}{n}\in\mathbb{Z}$. If $b>0$, then this follows
immediately from Lemma \ref{lem.sol.choose.a/b.sol2.pos}. Hence, for the rest
of this proof, we WLOG assume that we don't have $b>0$. Hence, $b\leq0$, so
that $b<0$ (since $b\neq0$). Therefore, $-b>0$. Thus, Lemma
\ref{lem.sol.choose.a/b.sol2.lem} (applied to $-a$ and $-b$ instead of $a$ and
$b$) yields $\left(  -b\right)  ^{2n-1}\dbinom{\left(  -a\right)  /\left(
-b\right)  }{n}\in\mathbb{Z}$. Since $\left(  -a\right)  /\left(  -b\right)
=a/b$, this rewrites as $\left(  -b\right)  ^{2n-1}\dbinom{a/b}{n}%
\in\mathbb{Z}$. In other words, $\left(  -b\right)  ^{2n-1}\dbinom{a/b}{n}$ is
an integer.

But%
\[
\underbrace{\left(  -b\right)  ^{2n-1}}_{=\left(  -1\right)  ^{2n-1}b^{2n-1}%
}\dbinom{a/b}{n}=\underbrace{\left(  -1\right)  ^{2n-1}}%
_{\substack{=-1\\\text{(since }2n-1\text{ is odd)}}}b^{2n-1}\dbinom{a/b}%
{n}=-b^{2n-1}\dbinom{a/b}{n}.
\]
In other words, the two numbers $\left(  -b\right)  ^{2n-1}\dbinom{a/b}{n}$
and $b^{2n-1}\dbinom{a/b}{n}$ differ only in sign. Since the first of them is
an integer, we thus conclude that so is the second. In other words,
$b^{2n-1}\dbinom{a/b}{n}\in\mathbb{Z}$. This proves Theorem
\ref{thm.sol.choose.a/b.sol2.main}.
\end{proof}
\end{vershort}

\begin{verlong}
\begin{proof}
[Proof of Theorem \ref{thm.sol.choose.a/b.sol2.main}.]We must prove that
$b^{2n-1}\dbinom{a/b}{n}\in\mathbb{Z}$. If $b>0$, then this follows
immediately from Lemma \ref{lem.sol.choose.a/b.sol2.pos}. Hence, for the rest
of this proof, we can WLOG assume that we don't have $b>0$. Assume this.

We have $b\leq0$ (since we don't have $b>0$). Combining this with $b\neq0$, we
obtain $b<0$. Hence, $-b>0$. Thus, Lemma \ref{lem.sol.choose.a/b.sol2.lem}
(applied to $-a$ and $-b$ instead of $a$ and $b$) yields $\left(  -b\right)
^{2n-1}\dbinom{\left(  -a\right)  /\left(  -b\right)  }{n}\in\mathbb{Z}$.
Since $\left(  -a\right)  /\left(  -b\right)  =a/b$, this rewrites as $\left(
-b\right)  ^{2n-1}\dbinom{a/b}{n}\in\mathbb{Z}$. In other words, $\left(
-b\right)  ^{2n-1}\dbinom{a/b}{n}$ is an integer.

Both numbers $\left(  -1\right)  ^{2n-1}$ and $\left(  -b\right)
^{2n-1}\dbinom{a/b}{n}$ are integers. Hence, their product is also an integer.
In other words, $\left(  -1\right)  ^{2n-1}\left(  -b\right)  ^{2n-1}%
\dbinom{a/b}{n}$ is an integer. In other words,
\[
\left(  -1\right)  ^{2n-1}\left(  -b\right)  ^{2n-1}\dbinom{a/b}{n}%
\in\mathbb{Z}.
\]

But $\left(  -1\right)  ^{2n-1}\left(  -b\right)  ^{2n-1}=\left(
\underbrace{\left(  -1\right)  \left(  -b\right)  }_{=b}\right)
^{2n-1}=b^{2n-1}$, so that \newline$b^{2n-1}=\left(  -1\right)  ^{2n-1}\left(
-b\right)  ^{2n-1}$. Multiplying this equality by $\dbinom{a/b}{n}$, we find%
\[
b^{2n-1}\dbinom{a/b}{n}=\left(  -1\right)  ^{2n-1}\left(  -b\right)
^{2n-1}\dbinom{a/b}{n}\in\mathbb{Z}.
\]
This proves Theorem \ref{thm.sol.choose.a/b.sol2.main}.
\end{proof}
\end{verlong}

Now, only some trivial bookkeeping remains to be done in order to solve
Exercise \ref{exe.choose.a/b}. To simplify it, we state the following
corollary from Theorem \ref{thm.sol.choose.a/b.sol2.main}:

\begin{corollary}
\label{cor.sol.choose.a/b.sol2.main}Let $a$ and $b$ be two integers such that
$b\neq0$. Let $n\in\mathbb{N}$. Let $m=\max\left\{  0,2n-1\right\}  $. Then,
$b^{m}\dbinom{a/b}{n}\in\mathbb{Z}$.
\end{corollary}

\begin{vershort}
\begin{proof}
[Proof of Corollary \ref{cor.sol.choose.a/b.sol2.main}.]If $n=0$, then
Corollary \ref{cor.sol.choose.a/b.sol2.main} holds\footnote{\textit{Proof.}
Assume that $n=0$. We must show that Corollary
\ref{cor.sol.choose.a/b.sol2.main} holds.
\par
We have $m=\max\left\{  0,2n-1\right\}  =0$ (since $2\underbrace{n}%
_{=0}-1=0-1<0$). Hence, $b^{m}=b^{0}=1$ and thus
\[
\underbrace{b^{m}}_{=1}\dbinom{a/b}{n}=\dbinom{a/b}{n}=\dbinom{a/b}%
{0}\ \ \ \ \ \ \ \ \ \ \left(  \text{since }n=0\right)  .
\]
\par
But Proposition \ref{prop.binom.00} \textbf{(a)} (applied to $a/b$ instead of
$m$) yields $\dbinom{a/b}{0}=1$. Thus, $b^{m}\dbinom{a/b}{n}=\dbinom{a/b}%
{0}=1\in\mathbb{Z}$. Thus, Corollary \ref{cor.sol.choose.a/b.sol2.main}
holds.}. Hence, for the rest of this proof, we WLOG assume that $n\neq0$.
Therefore, $n$ is a positive integer (since $n\in\mathbb{N}$). Therefore,
$2n-1>0$, so that $\max\left\{  0,2n-1\right\}  =2n-1$. Now, $m=\max\left\{
0,2n-1\right\}  =2n-1$, so that $b^{m}=b^{2n-1}$. Multiplying this equality by
$\dbinom{a/b}{n}$, we find%
\[
b^{m}\dbinom{a/b}{n}=b^{2n-1}\dbinom{a/b}{n}\in\mathbb{Z}%
\]
(by Theorem \ref{thm.sol.choose.a/b.sol2.main}). Corollary
\ref{cor.sol.choose.a/b.sol2.main} is thus proven.
\end{proof}
\end{vershort}

\begin{verlong}
\begin{proof}
[Proof of Corollary \ref{cor.sol.choose.a/b.sol2.main}.]If $n=0$, then
Corollary \ref{cor.sol.choose.a/b.sol2.main} holds\footnote{\textit{Proof.}
Assume that $n=0$. We must show that Corollary
\ref{cor.sol.choose.a/b.sol2.main} holds.
\par
We have $2\underbrace{n}_{=0}-1=0-1=-1$ and $m=\max\left\{
0,\underbrace{2n-1}_{=-1}\right\}  =\max\left\{  0,-1\right\}  =0$. Hence,
$b^{m}=b^{0}=1$ and thus $\underbrace{b^{m}}_{=1}\dbinom{a/b}{n}=\dbinom
{a/b}{n}=\dbinom{a/b}{0}$ (since $n=0$).
\par
But Proposition \ref{prop.binom.00} \textbf{(a)} (applied to $a/b$ instead of
$m$) yields $\dbinom{a/b}{0}=1$. Thus, $b^{m}\dbinom{a/b}{n}=\dbinom{a/b}%
{0}=1\in\mathbb{Z}$. Thus, Corollary \ref{cor.sol.choose.a/b.sol2.main} holds.
Qed.}. Hence, for the rest of this proof, we can WLOG assume that we don't
have $n=0$. Assume this.

We don't have $n=0$. Hence, we have $n\neq0$. Thus, $n$ is a positive integer
(since $n\in\mathbb{N}$). Therefore, $n\geq1$, so that $2\underbrace{n}%
_{\geq1}-1\geq2\cdot1-1=1>0$. Thus, $\max\left\{  0,2n-1\right\}  =2n-1$. Now,
$m=\max\left\{  0,2n-1\right\}  =2n-1$, so that $b^{m}=b^{2n-1}$. Multiplying
this equality by $\dbinom{a/b}{n}$, we find%
\[
b^{m}\dbinom{a/b}{n}=b^{2n-1}\dbinom{a/b}{n}\in\mathbb{Z}%
\]
(by Theorem \ref{thm.sol.choose.a/b.sol2.main}). Corollary
\ref{cor.sol.choose.a/b.sol2.main} is thus proven.
\end{proof}
\end{verlong}

\begin{proof}
[Second solution to Exercise \ref{exe.choose.a/b}.]Let $m=\max\left\{
0,2n-1\right\}  $. Then, \newline$m=\max\left\{  0,2n-1\right\}  \geq0$, so
that $m\in\mathbb{N}$. Also, Corollary \ref{cor.sol.choose.a/b.sol2.main}
shows that $b^{m}\dbinom{a/b}{n}\in\mathbb{Z}$. Hence, there exists some
$N\in\mathbb{N}$ such that $b^{N}\dbinom{a/b}{n}\in\mathbb{Z}$ (namely,
$N=m$). This solves Exercise \ref{exe.choose.a/b}.
\end{proof}

We remark that Theorem \ref{thm.sol.choose.a/b.sol2.main} is equivalent to
\cite[\S 3.3, problem 4]{AndDosS} (which is
\href{https://mks.mff.cuni.cz/kalva/short/soln/sh851.html}{a problem from the
IMO Shortlist 1985}). Indeed, the latter problem claims the following:

\begin{theorem}
\label{thm.sol.choose.a/b.sol2.ISL}Let $a$ and $b$ be integers. Let $n$ be a
positive integer. Then,%
\[
\dfrac{1}{n!}\cdot a\left(  a+b\right)  \left(  a+2b\right)  \cdots\left(
a+\left(  n-1\right)  b\right)  \cdot b^{n-1}\in\mathbb{Z}.
\]

\end{theorem}

Let us derive Theorem \ref{thm.sol.choose.a/b.sol2.ISL} from Theorem
\ref{thm.sol.choose.a/b.sol2.main}:

\begin{proof}
[Proof of Theorem \ref{thm.sol.choose.a/b.sol2.ISL}.]We are in one of the
following three cases:

\textit{Case 1:} We have $b\neq0$.

\textit{Case 2:} We have $n=1$.

\textit{Case 3:} We have neither $b\neq0$ nor $n=1$.

Let us first consider Case 1. In this case, we have $b\neq0$. Thus, $-b\neq0$
as well. Hence, Theorem \ref{thm.sol.choose.a/b.sol2.main} (applied to $-b$
instead of $b$) yields $\left(  -b\right)  ^{2n-1}\dbinom{a/\left(  -b\right)
}{n}\in\mathbb{Z}$.

But (\ref{eq.binom.mn}) (applied to $m=a/\left(  -b\right)  $) yields%
\begin{align*}
\dbinom{a/\left(  -b\right)  }{n}  &  =\dfrac{\left(  a/\left(  -b\right)
\right)  \left(  a/\left(  -b\right)  -1\right)  \cdots\left(  a/\left(
-b\right)  -n+1\right)  }{n!}\\
&  =\dfrac{1}{n!}\cdot\underbrace{\left(  a/\left(  -b\right)  \right)
\left(  a/\left(  -b\right)  -1\right)  \cdots\left(  a/\left(  -b\right)
-n+1\right)  }_{=\prod_{i=0}^{n-1}\left(  a/\left(  -b\right)  -i\right)  }\\
&  =\dfrac{1}{n!}\cdot\prod_{i=0}^{n-1}\underbrace{\left(  a/\left(
-b\right)  -i\right)  }_{=\dfrac{a+ib}{-b}}=\dfrac{1}{n!}\cdot
\underbrace{\prod_{i=0}^{n-1}\dfrac{a+ib}{-b}}_{=\dfrac{\prod_{i=0}%
^{n-1}\left(  a+ib\right)  }{\left(  -b\right)  ^{n}}}\\
&  =\dfrac{1}{n!}\cdot\dfrac{\prod_{i=0}^{n-1}\left(  a+ib\right)  }{\left(
-b\right)  ^{n}}.
\end{align*}
Multiplying both sides of this equality by $\left(  -b\right)  ^{2n-1}$, we
obtain%
\begin{align*}
\left(  -b\right)  ^{2n-1}\dbinom{a/\left(  -b\right)  }{n}  &  =\left(
-b\right)  ^{2n-1}\dfrac{1}{n!}\cdot\dfrac{\prod_{i=0}^{n-1}\left(
a+ib\right)  }{\left(  -b\right)  ^{n}}=\dfrac{1}{n!}\cdot\underbrace{\dfrac
{\left(  -b\right)  ^{2n-1}}{\left(  -b\right)  ^{n}}}_{=\left(  -b\right)
^{n-1}}\cdot\prod_{i=0}^{n-1}\left(  a+ib\right) \\
&  =\dfrac{1}{n!}\cdot\left(  -b\right)  ^{n-1}\cdot\prod_{i=0}^{n-1}\left(
a+ib\right)  .
\end{align*}
Multiplying both sides of this equality by $\left(  -1\right)  ^{n-1}$, we
obtain%
\begin{align*}
&  \left(  -1\right)  ^{n-1}\left(  -b\right)  ^{2n-1}\dbinom{a/\left(
-b\right)  }{n}\\
&  =\left(  -1\right)  ^{n-1}\dfrac{1}{n!}\cdot\left(  -b\right)  ^{n-1}%
\cdot\prod_{i=0}^{n-1}\left(  a+ib\right)  =\dfrac{1}{n!}\cdot
\underbrace{\left(  -1\right)  ^{n-1}\left(  -b\right)  ^{n-1}}_{=\left(
\left(  -1\right)  \left(  -b\right)  \right)  ^{n-1}}\cdot\prod_{i=0}%
^{n-1}\left(  a+ib\right) \\
&  =\dfrac{1}{n!}\cdot\left(  \underbrace{\left(  -1\right)  \left(
-b\right)  }_{=b}\right)  ^{n-1}\cdot\prod_{i=0}^{n-1}\left(  a+ib\right)
=\dfrac{1}{n!}\cdot b^{n-1}\cdot\prod_{i=0}^{n-1}\left(  a+ib\right) \\
&  =\dfrac{1}{n!}\cdot\underbrace{\prod_{i=0}^{n-1}\left(  a+ib\right)
}_{=a\left(  a+b\right)  \left(  a+2b\right)  \cdots\left(  a+\left(
n-1\right)  b\right)  }\cdot b^{n-1}\\
&  =\dfrac{1}{n!}\cdot a\left(  a+b\right)  \left(  a+2b\right)  \cdots\left(
a+\left(  n-1\right)  b\right)  \cdot b^{n-1}.
\end{align*}
Hence,%
\begin{align*}
&  \dfrac{1}{n!}\cdot a\left(  a+b\right)  \left(  a+2b\right)  \cdots\left(
a+\left(  n-1\right)  b\right)  \cdot b^{n-1}\\
&  =\underbrace{\left(  -1\right)  ^{n-1}}_{\in\mathbb{Z}}\underbrace{\left(
-b\right)  ^{2n-1}\dbinom{a/\left(  -b\right)  }{n}}_{\in\mathbb{Z}}%
\in\mathbb{Z}.
\end{align*}
Hence, Theorem \ref{thm.sol.choose.a/b.sol2.ISL} is proven in Case 1.

Let us now consider Case 2. In this case, we have $n=1$. Thus, $\dfrac{1}%
{n!}=\dfrac{1}{1!}=\dfrac{1}{1}=1\in\mathbb{Z}$ and $b^{n-1}=b^{1-1}%
=b^{0}=1\in\mathbb{Z}$, and therefore%
\[
\underbrace{\dfrac{1}{n!}}_{\in\mathbb{Z}}\cdot\underbrace{a\left(
a+b\right)  \left(  a+2b\right)  \cdots\left(  a+\left(  n-1\right)  b\right)
}_{\substack{\in\mathbb{Z}\\\text{(since }a\text{ and }b\text{ are integers)}%
}}\cdot\underbrace{b^{n-1}}_{\in\mathbb{Z}}\in\mathbb{Z}.
\]
Hence, Theorem \ref{thm.sol.choose.a/b.sol2.ISL} is proven in Case 2.

Let us finally consider Case 3. In this case, we have neither $b\neq0$ nor
$n=1$. Thus, $b=0$ (since we don't have $b\neq0$) and $n\neq1$ (since we don't
have $n=1$). From $n\neq1$, we obtain $n\geq2$ (since $n$ is a positive
integer), and thus $n-1\geq1$. Hence, $0^{n-1}=0$. But from $b=0$, we obtain
$b^{n-1}=0^{n-1}=0$. Thus,%
\[
\dfrac{1}{n!}\cdot a\left(  a+b\right)  \left(  a+2b\right)  \cdots\left(
a+\left(  n-1\right)  b\right)  \cdot\underbrace{b^{n-1}}_{=0}=0\in
\mathbb{Z}.
\]
Hence, Theorem \ref{thm.sol.choose.a/b.sol2.ISL} is proven in Case 3.

We have now proven Theorem \ref{thm.sol.choose.a/b.sol2.ISL} in each of the
three Cases 1, 2 and 3. Since these three Cases cover all possibilities, we
thus conclude that Theorem \ref{thm.sol.choose.a/b.sol2.ISL} always holds.
\end{proof}

Similarly, we can derive Theorem \ref{thm.sol.choose.a/b.sol2.main} back from
Theorem \ref{thm.sol.choose.a/b.sol2.ISL}. (That derivation is actually
easier, since we don't need to distinguish between any cases.)

\subsection{Solution to Exercise \ref{exe.ISL1975P7}}

We shall derive Exercise \ref{exe.ISL1975P7} from a sequence of lemmas. The
first one will be an easy consequence from the binomial identity:

\begin{lemma}
\label{lem.sol.ISL1975P7.binom}Let $g\in\mathbb{N}$. Let $x$ and $y$ be two
real numbers. Let $z=x+y$.

\textbf{(a)} Then,%
\[
z^{g}=\sum_{i=0}^{g}\dbinom{g}{i}x^{i}y^{g-i}.
\]

\textbf{(b)} Let $n\in\mathbb{N}$ be such that $n\geq g$. Then,%
\[
y^{n-g}z^{g}=\sum_{i=0}^{n}\dbinom{g}{i}x^{i}y^{n-i}.
\]

\end{lemma}

\begin{proof}
[Proof of Lemma \ref{lem.sol.ISL1975P7.binom}.]Proposition
\ref{prop.binom.binomial} (applied to $g$ instead of $n$) yields
\[
\left(  x+y\right)  ^{g}=\sum_{k=0}^{g}\dbinom{g}{k}x^{k}y^{g-k}=\sum
_{i=0}^{g}\dbinom{g}{i}x^{i}y^{g-i}%
\]
(here, we have renamed the summation index $k$ as $i$ in the sum). Now, from
$z=x+y$, we obtain%
\begin{equation}
z^{g}=\left(  x+y\right)  ^{g}=\sum_{i=0}^{g}\dbinom{g}{i}x^{i}y^{g-i}.
\label{pf.lem.sol.ISL1975P7.binom.a}%
\end{equation}
This proves Lemma \ref{lem.sol.ISL1975P7.binom} \textbf{(a)}.

\textbf{(b)} We have $n-g\geq0$ (since $n\geq g$), and thus $n-g\in\mathbb{N}%
$. Hence, $y^{n-g}$ is well-defined.

Multiplying both sides of the equality (\ref{pf.lem.sol.ISL1975P7.binom.a}) by
$y^{n-g}$, we obtain%
\begin{align}
y^{n-g}z^{g}  &  =y^{n-g}\sum_{i=0}^{g}\dbinom{g}{i}x^{i}y^{g-i}=\sum
_{i=0}^{g}\dbinom{g}{i}x^{i}\underbrace{y^{n-g}y^{g-i}}_{\substack{=y^{\left(
n-g\right)  +\left(  g-i\right)  }=y^{n-i}\\\text{(since }\left(  n-g\right)
+\left(  g-i\right)  =n-i\text{)}}}\nonumber\\
&  =\sum_{i=0}^{g}\dbinom{g}{i}x^{i}y^{n-i}.
\label{pf.lem.sol.ISL1975P7.binom.L}%
\end{align}

On the other hand, for each $i\in\left\{  g+1,g+2,\ldots,n\right\}  $, we have%
\begin{equation}
\dbinom{g}{i}=0 \label{pf.lem.sol.ISL1975P7.binom.0a}%
\end{equation}
\footnote{\textit{Proof of (\ref{pf.lem.sol.ISL1975P7.binom.0a}):} Let
$i\in\left\{  g+1,g+2,\ldots,n\right\}  $. Thus, $i\geq g+1>g$, so that $g<i$.
Also, $i>g\geq0$ (since $g\in\mathbb{N}$), so that $i\in\mathbb{N}$. Hence,
Proposition \ref{prop.binom.0} (applied to $g$ and $i$ instead of $m$ and $n$)
yields $\dbinom{g}{i}=0$. This proves (\ref{pf.lem.sol.ISL1975P7.binom.0a}).}.

But $g\geq0$ (since $g\in\mathbb{N}$), so that $0\leq g\leq n$ (since $n\geq
g$). Hence, we can split the sum $\sum_{i=0}^{n}\dbinom{g}{i}x^{i}y^{n-i}$ as
follows:%
\begin{align*}
\sum_{i=0}^{n}\dbinom{g}{i}x^{i}y^{n-i}  &  =\sum_{i=0}^{g}\dbinom{g}{i}%
x^{i}y^{n-i}+\sum_{i=g+1}^{n}\underbrace{\dbinom{g}{i}}%
_{\substack{=0\\\text{(by (\ref{pf.lem.sol.ISL1975P7.binom.0a}))}}%
}x^{i}y^{n-i}\\
&  =\sum_{i=0}^{g}\dbinom{g}{i}x^{i}y^{n-i}+\underbrace{\sum_{i=g+1}^{n}%
0x^{i}y^{n-i}}_{=0}=\sum_{i=0}^{g}\dbinom{g}{i}x^{i}y^{n-i}.
\end{align*}
Comparing this with (\ref{pf.lem.sol.ISL1975P7.binom.L}), we obtain
$y^{n-g}z^{g}=\sum_{i=0}^{n}\dbinom{g}{i}x^{i}y^{n-i}$. This proves Lemma
\ref{lem.sol.ISL1975P7.binom} \textbf{(b)}.
\end{proof}

Our next lemma is a combinatorial identity that follows easily from
Proposition \ref{prop.vandermonde.consequences} \textbf{(f)}:

\begin{lemma}
\label{lem.sol.ISL1975P7.vdm}Let $n\in\mathbb{N}$ and $m\in\mathbb{N}$. Let
$i\in\left\{  0,1,\ldots,n\right\}  $. Then,%
\[
\sum_{k=0}^{n}\dbinom{m+k}{k}\dbinom{n-k}{i}=\dbinom{n+m+1}{m+i+1}.
\]

\end{lemma}

\begin{proof}
[Proof of Lemma \ref{lem.sol.ISL1975P7.vdm}.]We have $i\in\left\{
0,1,\ldots,n\right\}  \subseteq\mathbb{N}$. Also, from $i\in\left\{
0,1,\ldots,n\right\}  $, we obtain $i\leq n$. Hence, $m+\underbrace{i}_{\leq
n}+1\leq m+n+1=n+m+1$. In other words, $n+m+1\geq m+i+1$. Also, $m+i+1\in
\mathbb{N}$ (since $m\in\mathbb{N}$ and $i\in\mathbb{N}$) and $n+m+1\in
\mathbb{N}$ (since $n\in\mathbb{N}$ and $m\in\mathbb{N}$). Thus, Proposition
\ref{prop.binom.symm} (applied to $n+m+1$ and $m+i+1$ instead of $m$ and $n$)
yields%
\begin{equation}
\dbinom{n+m+1}{m+i+1}=\dbinom{n+m+1}{\left(  n+m+1\right)  -\left(
m+i+1\right)  }=\dbinom{n+m+1}{n-i} \label{pf.lem.sol.ISL1975P7.vdm.sym1}%
\end{equation}
(since $\left(  n+m+1\right)  -\left(  m+i+1\right)  =n-i$).

For each $k\in\mathbb{N}$, we have%
\begin{equation}
\dbinom{m+k}{m}=\dbinom{m+k}{k} \label{pf.lem.sol.ISL1975P7.vdm.sym2}%
\end{equation}
(by Lemma \ref{lem.binom.symmetry-m+n} (applied to $k$ instead of $n$)).

For each $k\in\left\{  0,1,\ldots,m-1\right\}  $, we have%
\begin{equation}
\dbinom{k}{m}=0 \label{pf.lem.sol.ISL1975P7.vdm.3}%
\end{equation}
\footnote{\textit{Proof of (\ref{pf.lem.sol.ISL1975P7.vdm.3}):} Let
$k\in\left\{  0,1,\ldots,m-1\right\}  $. Thus, $k\geq0$ and $k\leq m-1<m$.
From $k\geq0$, we obtain $k\in\mathbb{N}$. Hence, Proposition
\ref{prop.binom.0} (applied to $k$ and $m$ instead of $m$ and $n$) yields
$\dbinom{k}{m}=0$. This proves (\ref{pf.lem.sol.ISL1975P7.vdm.3}).}.

We have $m\geq0$ (since $m\in\mathbb{N}$) and $n\geq0$ (since $n\in\mathbb{N}%
$), so that $\underbrace{n}_{\geq0}+m\geq m$ and thus $m\leq n+m$. From
$m\geq0$, we obtain $0\leq m\leq n+m$. Hence, we can split the sum $\sum
_{k=0}^{n+m}\dbinom{k}{m}\dbinom{n+m-k}{i}$ as follows:%
\begin{align}
&  \sum_{k=0}^{n+m}\dbinom{k}{m}\dbinom{n+m-k}{i}\nonumber\\
&  =\sum_{k=0}^{m-1}\underbrace{\dbinom{k}{m}}_{\substack{=0\\\text{(by
(\ref{pf.lem.sol.ISL1975P7.vdm.3}))}}}\dbinom{n+m-k}{i}+\sum_{k=m}%
^{n+m}\dbinom{k}{m}\dbinom{n+m-k}{i}\nonumber\\
&  =\underbrace{\sum_{k=0}^{m-1}0\dbinom{n+m-k}{i}}_{=0}+\sum_{k=m}%
^{n+m}\dbinom{k}{m}\dbinom{n+m-k}{i}\nonumber\\
&  =\sum_{k=m}^{n+m}\dbinom{k}{m}\dbinom{n+m-k}{i}=\sum_{k=0}^{n}%
\underbrace{\dbinom{k+m}{m}}_{\substack{=\dbinom{m+k}{m}\\\text{(since
}k+m=m+k\text{)}}}\ \ \underbrace{\dbinom{n+m-\left(  k+m\right)  }{i}%
}_{\substack{=\dbinom{n-k}{i}\\\text{(since }n+m-\left(  k+m\right)
=n-k\text{)}}}\nonumber\\
&  \ \ \ \ \ \ \ \ \ \ \left(  \text{here, we have substituted }k+m\text{ for
}k\text{ in the sum}\right) \nonumber\\
&  =\sum_{k=0}^{n}\underbrace{\dbinom{m+k}{m}}_{\substack{=\dbinom{m+k}%
{k}\\\text{(by (\ref{pf.lem.sol.ISL1975P7.vdm.sym2}))}}}\dbinom{n-k}{i}%
=\sum_{k=0}^{n}\dbinom{m+k}{k}\dbinom{n-k}{i}.
\label{pf.lem.sol.ISL1975P7.vdm.5}%
\end{align}

But Proposition \ref{prop.vandermonde.consequences} \textbf{(f)} (applied to
$n+m$, $m$ and $i$ instead of $n$, $x$ and $y$) yields%
\[
\dbinom{n+m+1}{m+i+1}=\sum_{k=0}^{n+m}\dbinom{k}{m}\dbinom{n+m-k}{i}%
=\sum_{k=0}^{n}\dbinom{m+k}{k}\dbinom{n-k}{i}%
\]
(by (\ref{pf.lem.sol.ISL1975P7.vdm.5})). This proves Lemma
\ref{lem.sol.ISL1975P7.vdm}.
\end{proof}

Our next step brings us much closer to Exercise \ref{exe.ISL1975P7}:

\begin{proposition}
\label{prop.sol.ISL1975P7.gen}Let $x$ and $y$ be two real numbers. Let
$z=x+y$. Let $n\in\mathbb{N}$ and $m\in\mathbb{N}$. Then:

\textbf{(a)} We have%
\[
x^{m+1}\sum_{k=0}^{n}\dbinom{m+k}{k}y^{k}z^{n-k}=\sum_{i=m+1}^{n+m+1}%
\dbinom{n+m+1}{i}x^{i}y^{\left(  n+m+1\right)  -i}.
\]

\textbf{(b)} We have%
\[
y^{n+1}\sum_{k=0}^{m}\dbinom{n+k}{k}x^{k}z^{m-k}=\sum_{i=0}^{m}\dbinom
{n+m+1}{i}x^{i}y^{\left(  n+m+1\right)  -i}.
\]

\textbf{(c)} We have%
\[
x^{m+1}\sum_{k=0}^{n}\dbinom{m+k}{k}y^{k}z^{n-k}+y^{n+1}\sum_{k=0}^{m}%
\dbinom{n+k}{k}x^{k}z^{m-k}=z^{n+m+1}.
\]

\end{proposition}

\begin{proof}
[Proof of Proposition \ref{prop.sol.ISL1975P7.gen}.]\textbf{(a)} Let
$k\in\left\{  0,1,\ldots,n\right\}  $. Thus, $k\leq n$, so that $n-k\geq0$ and
thus $n-k\in\mathbb{N}$. Also, $k\geq0$ (since $k\in\left\{  0,1,\ldots
,n\right\}  $), so that $n-\underbrace{k}_{\geq0}\leq n$ and thus $n\geq n-k$.
Hence, Lemma \ref{lem.sol.ISL1975P7.binom} \textbf{(a)} (applied to $g=n-k$)
yields%
\[
y^{n-\left(  n-k\right)  }z^{n-k}=\sum_{i=0}^{n}\dbinom{n-k}{i}x^{i}y^{n-i}.
\]
In view of $n-\left(  n-k\right)  =k$, this rewrites as follows:%
\begin{equation}
y^{k}z^{n-k}=\sum_{i=0}^{n}\dbinom{n-k}{i}x^{i}y^{n-i}.
\label{pf.prop.sol.ISL1975P7.gen.a.1}%
\end{equation}

Now, forget that we fixed $k$. We thus have proven the equality
(\ref{pf.prop.sol.ISL1975P7.gen.a.1}) for each $k\in\left\{  0,1,\ldots
,n\right\}  $. Now,%
\begin{align*}
&  x^{m+1}\sum_{k=0}^{n}\dbinom{m+k}{k}\underbrace{y^{k}z^{n-k}}%
_{\substack{=\sum_{i=0}^{n}\dbinom{n-k}{i}x^{i}y^{n-i}\\\text{(by
(\ref{pf.prop.sol.ISL1975P7.gen.a.1}))}}}\\
&  =x^{m+1}\sum_{k=0}^{n}\underbrace{\dbinom{m+k}{k}\sum_{i=0}^{n}\dbinom
{n-k}{i}x^{i}y^{n-i}}_{=\sum_{i=0}^{n}\dbinom{m+k}{k}\dbinom{n-k}{i}%
x^{i}y^{n-i}}=x^{m+1}\underbrace{\sum_{k=0}^{n}\sum_{i=0}^{n}}_{=\sum
_{i=0}^{n}\sum_{k=0}^{n}}\dbinom{m+k}{k}\dbinom{n-k}{i}x^{i}y^{n-i}\\
&  =x^{m+1}\sum_{i=0}^{n}\underbrace{\sum_{k=0}^{n}\dbinom{m+k}{k}\dbinom
{n-k}{i}x^{i}y^{n-i}}_{=\left(  \sum_{k=0}^{n}\dbinom{m+k}{k}\dbinom{n-k}%
{i}\right)  x^{i}y^{n-i}}=x^{m+1}\sum_{i=0}^{n}\underbrace{\left(  \sum
_{k=0}^{n}\dbinom{m+k}{k}\dbinom{n-k}{i}\right)  }_{\substack{=\dbinom
{n+m+1}{m+i+1}\\\text{(by Lemma \ref{lem.sol.ISL1975P7.vdm})}}}x^{i}y^{n-i}\\
&  =x^{m+1}\sum_{i=0}^{n}\underbrace{\dbinom{n+m+1}{m+i+1}}%
_{\substack{=\dbinom{n+m+1}{\left(  m+1\right)  +i}\\\text{(since
}m+i+1=\left(  m+1\right)  +i\text{)}}}x^{i}\underbrace{y^{n-i}}%
_{\substack{=y^{\left(  n+m+1\right)  -\left(  \left(  m+1\right)  +i\right)
}\\\text{(since }n-i=\left(  n+m+1\right)  -\left(  \left(  m+1\right)
+i\right)  \text{)}}}\\
&  =x^{m+1}\sum_{i=0}^{n}\dbinom{n+m+1}{\left(  m+1\right)  +i}x^{i}y^{\left(
n+m+1\right)  -\left(  \left(  m+1\right)  +i\right)  }\\
&  =\sum_{i=0}^{n}\dbinom{n+m+1}{\left(  m+1\right)  +i}\underbrace{x^{m+1}%
x^{i}}_{=x^{\left(  m+1\right)  +i}}y^{\left(  n+m+1\right)  -\left(  \left(
m+1\right)  +i\right)  }\\
&  =\sum_{i=0}^{n}\dbinom{n+m+1}{\left(  m+1\right)  +i}x^{\left(  m+1\right)
+i}y^{\left(  n+m+1\right)  -\left(  \left(  m+1\right)  +i\right)  }\\
&  =\sum_{i=m+1}^{\left(  m+1\right)  +n}\dbinom{n+m+1}{i}x^{i}y^{\left(
n+m+1\right)  -i}\\
&  \ \ \ \ \ \ \ \ \ \ \left(  \text{here, we have substituted }i\text{ for
}\left(  m+1\right)  +i\text{ in the sum}\right) \\
&  =\sum_{i=m+1}^{n+m+1}\dbinom{n+m+1}{i}x^{i}y^{\left(  n+m+1\right)
-i}\ \ \ \ \ \ \ \ \ \ \left(  \text{since }\left(  m+1\right)
+n=n+m+1\right)  .
\end{align*}
This proves Proposition \ref{prop.sol.ISL1975P7.gen} \textbf{(a)}.

\textbf{(b)} Each $i\in\left\{  n+1,n+2,\ldots,n+m+1\right\}  $ satisfies%
\begin{equation}
\dbinom{n+m+1}{i}=\dbinom{n+m+1}{\left(  n+m+1\right)  -i}
\label{pf.prop.sol.ISL1975P7.gen.b.1}%
\end{equation}
\footnote{\textit{Proof of (\ref{pf.prop.sol.ISL1975P7.gen.b.1}):} Let
$i\in\left\{  n+1,n+2,\ldots,n+m+1\right\}  $. Thus, $i\geq n+1>n\geq0$ (since
$n\in\mathbb{N}$). Hence, $i\in\mathbb{N}$. Also, $i\in\left\{  n+1,n+2,\ldots
,n+m+1\right\}  $ shows that $i\leq n+m+1$, so that $n+m+1\geq i$. Thus,
Proposition \ref{prop.binom.symm} (applied to $n+m+1$ and $i$ instead of $m$
and $n$) yields $\dbinom{n+m+1}{i}=\dbinom{n+m+1}{\left(  n+m+1\right)  -i}$.
This proves (\ref{pf.prop.sol.ISL1975P7.gen.b.1}).}.

We have $z=x+y=y+x$. Hence, Proposition \ref{prop.sol.ISL1975P7.gen}
\textbf{(a)} (applied to $y$, $x$, $m$ and $n$ instead of $x$, $y$, $n$ and
$m$) yields%
\begin{align*}
&  y^{n+1}\sum_{k=0}^{m}\dbinom{n+k}{k}x^{k}z^{m-k}\\
&  =\sum_{i=n+1}^{m+n+1}\dbinom{m+n+1}{i}y^{i}x^{\left(  m+n+1\right)  -i}\\
&  =\sum_{i=n+1}^{n+m+1}\underbrace{\dbinom{n+m+1}{i}}_{\substack{=\dbinom
{n+m+1}{\left(  n+m+1\right)  -i}\\\text{(by
(\ref{pf.prop.sol.ISL1975P7.gen.b.1}))}}}\ \ \underbrace{y^{i}}%
_{\substack{=y^{\left(  n+m+1\right)  -\left(  \left(  n+m+1\right)
-i\right)  }\\\text{(since }i=\left(  n+m+1\right)  -\left(  \left(
n+m+1\right)  -i\right)  \text{)}}}x^{\left(  n+m+1\right)  -i}\\
&  \ \ \ \ \ \ \ \ \ \ \left(  \text{since }m+n=n+m\right) \\
&  =\sum_{i=n+1}^{n+m+1}\dbinom{n+m+1}{\left(  n+m+1\right)  -i}y^{\left(
n+m+1\right)  -\left(  \left(  n+m+1\right)  -i\right)  }x^{\left(
n+m+1\right)  -i}\\
&  =\underbrace{\sum_{i=0}^{\left(  n+m+1\right)  -\left(  n+1\right)  }%
}_{\substack{=\sum_{i=0}^{m}\\\text{(since }\left(  n+m+1\right)  -\left(
n+1\right)  =m\text{)}}}\dbinom{n+m+1}{i}\underbrace{y^{\left(  n+m+1\right)
-i}x^{i}}_{=x^{i}y^{\left(  n+m+1\right)  -i}}\\
&  \ \ \ \ \ \ \ \ \ \ \left(  \text{here, we have substituted }i\text{ for
}\left(  n+m+1\right)  -i\text{ in the sum}\right) \\
&  =\sum_{i=0}^{m}\dbinom{n+m+1}{i}x^{i}y^{\left(  n+m+1\right)  -i}.
\end{align*}
This proves Proposition \ref{prop.sol.ISL1975P7.gen} \textbf{(b)}.

\textbf{(c)} We have $\underbrace{n}_{\substack{\geq0\\\text{(since }%
n\in\mathbb{N}\text{)}}}+m+1\geq m+1\geq m$ and thus $m\leq n+m+1$. Also, from
$m\in\mathbb{N}$, we obtain $m\geq0$, thus $0\leq m\leq n+m+1$.

We have $n+m+1\in\mathbb{N}$ (since $n\in\mathbb{N}$ and $m\in\mathbb{N}$).
Hence, Lemma \ref{lem.sol.ISL1975P7.binom} \textbf{(a)} (applied to $n+m+1$
instead of $g$) yields%
\begin{align}
&  z^{n+m+1}\nonumber\\
&  =\sum_{i=0}^{n+m+1}\dbinom{n+m+1}{i}x^{i}y^{\left(  n+m+1\right)
-i}\nonumber\\
&  =\sum_{i=0}^{m}\dbinom{n+m+1}{i}x^{i}y^{\left(  n+m+1\right)  -i}%
+\sum_{i=m+1}^{n+m+1}\dbinom{n+m+1}{i}x^{i}y^{\left(  n+m+1\right)  -i}
\label{pf.prop.sol.ISL1975P7.gen.b.split}%
\end{align}
(here, we have split the sum at $i=m$, since $0\leq m\leq n+m+1$).

Proposition \ref{prop.sol.ISL1975P7.gen} \textbf{(a)} yields%
\begin{equation}
x^{m+1}\sum_{k=0}^{n}\dbinom{m+k}{k}y^{k}z^{n-k}=\sum_{i=m+1}^{n+m+1}%
\dbinom{n+m+1}{i}x^{i}y^{\left(  n+m+1\right)  -i}.
\label{pf.prop.sol.ISL1975P7.gen.b.a}%
\end{equation}
Proposition \ref{prop.sol.ISL1975P7.gen} \textbf{(b)} yields%
\[
y^{n+1}\sum_{k=0}^{m}\dbinom{n+k}{k}x^{k}z^{m-k}=\sum_{i=0}^{m}\dbinom
{n+m+1}{i}x^{i}y^{\left(  n+m+1\right)  -i}.
\]
Adding this equality to (\ref{pf.prop.sol.ISL1975P7.gen.b.a}), we find%
\begin{align*}
&  x^{m+1}\sum_{k=0}^{n}\dbinom{m+k}{k}y^{k}z^{n-k}+y^{n+1}\sum_{k=0}%
^{m}\dbinom{n+k}{k}x^{k}z^{m-k}\\
&  =\sum_{i=m+1}^{n+m+1}\dbinom{n+m+1}{i}x^{i}y^{\left(  n+m+1\right)
-i}+\sum_{i=0}^{m}\dbinom{n+m+1}{i}x^{i}y^{\left(  n+m+1\right)  -i}\\
&  =\sum_{i=0}^{m}\dbinom{n+m+1}{i}x^{i}y^{\left(  n+m+1\right)  -i}%
+\sum_{i=m+1}^{n+m+1}\dbinom{n+m+1}{i}x^{i}y^{\left(  n+m+1\right)  -i}\\
&  =z^{n+m+1}\ \ \ \ \ \ \ \ \ \ \left(  \text{by
(\ref{pf.prop.sol.ISL1975P7.gen.b.split})}\right)  .
\end{align*}
This proves Proposition \ref{prop.sol.ISL1975P7.gen} \textbf{(c)}.
\end{proof}

\begin{proof}
[Solution to Exercise \ref{exe.ISL1975P7}.]\textbf{(a)} Let $x$ and $y$ be two
real numbers such that $x+y=1$. Let $n\in\mathbb{N}$ and $m\in\mathbb{N}$.
Thus, $1=x+y$. Hence, Proposition \ref{prop.sol.ISL1975P7.gen} \textbf{(c)}
(applied to $z=1$) yields%
\[
x^{m+1}\sum_{k=0}^{n}\dbinom{m+k}{k}y^{k}1^{n-k}+y^{n+1}\sum_{k=0}^{m}%
\dbinom{n+k}{k}x^{k}1^{m-k}=1^{n+m+1}=1.
\]
Comparing this with%
\begin{align*}
&  x^{m+1}\sum_{k=0}^{n}\dbinom{m+k}{k}y^{k}\underbrace{1^{n-k}}_{=1}%
+y^{n+1}\sum_{k=0}^{m}\dbinom{n+k}{k}x^{k}\underbrace{1^{m-k}}_{=1}\\
&  =x^{m+1}\sum_{k=0}^{n}\dbinom{m+k}{k}y^{k}+y^{n+1}\sum_{k=0}^{m}%
\dbinom{n+k}{k}x^{k},
\end{align*}
we obtain%
\[
x^{m+1}\sum_{k=0}^{n}\dbinom{m+k}{k}y^{k}+y^{n+1}\sum_{k=0}^{m}\dbinom{n+k}%
{k}x^{k}=1.
\]
This solves Exercise \ref{exe.ISL1975P7} \textbf{(a)}.

\textbf{(b)} The real numbers $\dfrac{1}{2}$ and $\dfrac{1}{2}$ satisfy
$\dfrac{1}{2}+\dfrac{1}{2}=1$. Hence, Exercise \ref{exe.ISL1975P7}
\textbf{(a)} (applied to $m=n$, $x=\dfrac{1}{2}$ and $y=\dfrac{1}{2}$) yields%
\[
\left(  \dfrac{1}{2}\right)  ^{n+1}\sum_{k=0}^{n}\dbinom{n+k}{k}\left(
\dfrac{1}{2}\right)  ^{k}+\left(  \dfrac{1}{2}\right)  ^{n+1}\sum_{k=0}%
^{n}\dbinom{n+k}{k}\left(  \dfrac{1}{2}\right)  ^{k}=1.
\]
Hence,%
\begin{align*}
1  &  =\left(  \dfrac{1}{2}\right)  ^{n+1}\sum_{k=0}^{n}\dbinom{n+k}{k}\left(
\dfrac{1}{2}\right)  ^{k}+\left(  \dfrac{1}{2}\right)  ^{n+1}\sum_{k=0}%
^{n}\dbinom{n+k}{k}\left(  \dfrac{1}{2}\right)  ^{k}\\
&  =2\cdot\underbrace{\left(  \dfrac{1}{2}\right)  ^{n+1}}_{=\dfrac{1}%
{2^{n+1}}}\sum_{k=0}^{n}\dbinom{n+k}{k}\underbrace{\left(  \dfrac{1}%
{2}\right)  ^{k}}_{=\dfrac{1}{2^{k}}}=\underbrace{2\cdot\dfrac{1}{2^{n+1}}%
}_{=\dfrac{1}{2^{n}}}\sum_{k=0}^{n}\dbinom{n+k}{k}\dfrac{1}{2^{k}}=\dfrac
{1}{2^{n}}\sum_{k=0}^{n}\dbinom{n+k}{k}\dfrac{1}{2^{k}}.
\end{align*}
Multiplying both sides of this equality by $2^{n}$, we obtain $2^{n}%
=\sum_{k=0}^{n}\dbinom{n+k}{k}\dfrac{1}{2^{k}}$. This solves Exercise
\ref{exe.ISL1975P7} \textbf{(b)}.
\end{proof}

\subsection{Solution to Exercise \ref{exe.ps2.2.1}}

\begin{proof}
[Solution to Exercise \ref{exe.ps2.2.1}.]We claim that every $n\in\mathbb{N}$
satisfies%
\begin{equation}
x_{n}=\dfrac{1}{2^{n}}\left(  2na^{n-1}x_{1}-\left(  n-1\right)  a^{n}%
x_{0}\right)  \label{sol.ps2.2.1.claim}%
\end{equation}
(where $na^{n-1}$ is to be understood as $0$ when $n=0$\ \ \ \ \footnote{This
needs to be said, because $a^{n-1}$ alone can be undefined for $n=0$ (if
$a=0$).}).

[\textit{Proof of (\ref{sol.ps2.2.1.claim}):} We shall prove
(\ref{sol.ps2.2.1.claim}) by strong induction\footnote{See Section
\ref{sect.ind.SIP} for an introduction to strong induction.} on $n$. So we fix
some $N\in\mathbb{N}$, and we assume that (\ref{sol.ps2.2.1.claim}) is already
proven for every $n<N$. (This is our induction hypothesis.) We now need to
show that (\ref{sol.ps2.2.1.claim}) holds for $n=N$ as well.\footnote{If you
are wondering \textquotedblleft where is the induction base?\textquotedblright%
: It isn't missing. A strong induction needs no induction base (see Convention
\ref{conv.ind.SIPlang} for the details). Strong induction lets you prove that
some statement $\mathcal{A}_{n}$ holds for every $n\in\mathbb{N}$ by means of
proving that for every $N\in\mathbb{N}$,%
\begin{equation}
\text{if }\mathcal{A}_{n}\text{ holds for every }n<N\text{, then }%
\mathcal{A}_{N}\text{ holds.} \label{sol.ps2.2.1.strind}%
\end{equation}
This immediately shows that $\mathcal{A}_{0}$ holds: Namely, it is clear that
$\mathcal{A}_{n}$ holds for every $n<0$ (because there exists no $n<0$), and
thus (\ref{sol.ps2.2.1.strind}) (applied to $N=0$) shows that $\mathcal{A}%
_{0}$ holds.
\par
Of course, the proof of (\ref{sol.ps2.2.1.strind}) might involve some case
analysis; in particular, it might argue differently depending on whether $N=0$
or $N\geq1$. (Indeed, our proof will be something like this: it will treat the
cases $N=0$, $N=1$ and $N\geq2$ separately.) So there can be a
\textquotedblleft de-facto induction base\textquotedblright\ (or two, or many)
hidden in the proof of (\ref{sol.ps2.2.1.strind}).} In other words, we need to
prove that
\begin{equation}
x_{N}=\dfrac{1}{2^{N}}\left(  2Na^{N-1}x_{1}-\left(  N-1\right)  a^{N}%
x_{0}\right)  . \label{sol.ps2.2.1.goal}%
\end{equation}

We must be in one of the following three cases:

\textit{Case 1:} We have $N=0$.

\textit{Case 2:} We have $N=1$.

\textit{Case 3:} We have $N\geq2$.

Let us first consider Case 1. In this case, we have $N=0$. Hence,%
\begin{align*}
\dfrac{1}{2^{N}}\left(  2Na^{N-1}x_{1}-\left(  N-1\right)  a^{N}x_{0}\right)
&  =\underbrace{\dfrac{1}{2^{0}}}_{=1}\left(  \underbrace{2\cdot0a^{0-1}x_{1}%
}_{=0}-\underbrace{\left(  0-1\right)  a^{0}}_{=-1}x_{0}\right) \\
&  =1\left(  0-\left(  -1\right)  x_{0}\right)  =1x_{0}=x_{0}\\
&  =x_{N}\ \ \ \ \ \ \ \ \ \ \left(  \text{since }0=N\right)  .
\end{align*}
In other words, $x_{N}=\dfrac{1}{2^{N}}\left(  2Na^{N-1}x_{1}-\left(
N-1\right)  a^{N}x_{0}\right)  $. Hence, (\ref{sol.ps2.2.1.goal}) is proven in
Case 1.

The proof of (\ref{sol.ps2.2.1.goal}) in Case 2 is similarly straightforward,
and is left to the reader.

Let us now consider Case 3. In this case, we have $N\geq2$. Hence, both $N-1$
and $N-2$ are nonnegative integers. Moreover, $N-1<N$, so that
(\ref{sol.ps2.2.1.claim}) is already proven for $n=N-1$ (by our induction
hypothesis). In other words, we have%
\begin{equation}
x_{N-1}=\dfrac{1}{2^{N-1}}\left(  2\left(  N-1\right)  a^{N-2}x_{1}-\left(
N-2\right)  a^{N-1}x_{0}\right)  . \label{sol.ps2.2.1.hyp1}%
\end{equation}
Also, $N-2$ is a nonnegative integer such that $N-2<N$. Hence,
(\ref{sol.ps2.2.1.claim}) is already proven for $n=N-2$ (by our induction
hypothesis). In other words, we have%
\begin{equation}
x_{N-2}=\dfrac{1}{2^{N-2}}\left(  2\left(  N-2\right)  a^{N-3}x_{1}-\left(
N-3\right)  a^{N-2}x_{0}\right)  . \label{sol.ps2.2.1.hyp2}%
\end{equation}

Now, recall that $a^{2}+4b=0$, so that $4b=-a^{2}$. Hence, $4b\cdot\left(
N-2\right)  a^{N-3}=\left(  -a^{2}\right)  \cdot\left(  N-2\right)
a^{N-3}=-\left(  N-2\right)  a^{N-1}$. (Don't forget to check that this latter
equality holds also when $N-2=0$; keep in mind that $\left(  N-2\right)
a^{N-3}$ was defined to be $0$ in this case, although $a^{N-3}$ might be
undefined.) From $4b=-a^{2}$, we also deduce $b=-\dfrac{a^{2}}{4}$, so that
$b\left(  N-3\right)  a^{N-2}=\dfrac{-a^{2}}{4}\left(  N-3\right)
a^{N-2}=-\dfrac{1}{4}\left(  N-3\right)  a^{N}$.

But the sequence $\left(  x_{0},x_{1},x_{2},\ldots\right)  $ is $\left(
a,b\right)  $-recurrent. Hence,%
\begin{align*}
x_{N}  &  =ax_{N-1}+bx_{N-2}\\
&  =a\cdot\dfrac{1}{2^{N-1}}\left(  2\left(  N-1\right)  a^{N-2}x_{1}-\left(
N-2\right)  a^{N-1}x_{0}\right) \\
&  \ \ \ \ \ \ \ \ \ \ +b\cdot\dfrac{1}{2^{N-2}}\left(  2\left(  N-2\right)
a^{N-3}x_{1}-\left(  N-3\right)  a^{N-2}x_{0}\right) \\
&  \ \ \ \ \ \ \ \ \ \ \left(  \text{by (\ref{sol.ps2.2.1.hyp1}) and
(\ref{sol.ps2.2.1.hyp2})}\right) \\
&  =a\cdot\dfrac{1}{2^{N-1}}2\left(  N-1\right)  a^{N-2}x_{1}-a\cdot\dfrac
{1}{2^{N-1}}\left(  N-2\right)  a^{N-1}x_{0}\\
&  \ \ \ \ \ \ \ \ \ \ +b\cdot\dfrac{1}{2^{N-2}}2\left(  N-2\right)
a^{N-3}x_{1}-b\cdot\dfrac{1}{2^{N-2}}\left(  N-3\right)  a^{N-2}x_{0}\\
&  =\underbrace{\left(  a\cdot\dfrac{1}{2^{N-1}}2\left(  N-1\right)
a^{N-2}+b\cdot\dfrac{1}{2^{N-2}}2\left(  N-2\right)  a^{N-3}\right)
}_{=\dfrac{1}{2^{N-1}}\left(  a\cdot2\left(  N-1\right)  a^{N-2}%
+4b\cdot\left(  N-2\right)  a^{N-3}\right)  }x_{1}\\
&  \ \ \ \ \ \ \ \ \ \ -\underbrace{\left(  a\cdot\dfrac{1}{2^{N-1}}\left(
N-2\right)  a^{N-1}+b\cdot\dfrac{1}{2^{N-2}}\left(  N-3\right)  a^{N-2}%
\right)  }_{=\dfrac{1}{2^{N-1}}\left(  a\cdot\left(  N-2\right)
a^{N-1}+2b\left(  N-3\right)  a^{N-2}\right)  }x_{0}\\
&  =\dfrac{1}{2^{N-1}}\left(  \underbrace{a\cdot2\left(  N-1\right)  a^{N-2}%
}_{=2\left(  N-1\right)  a^{N-1}}+\underbrace{4b\cdot\left(  N-2\right)
a^{N-3}}_{=-\left(  N-2\right)  a^{N-1}}\right)  x_{1}\\
&  \ \ \ \ \ \ \ \ \ \ -\dfrac{1}{2^{N-1}}\left(  \underbrace{a\cdot\left(
N-2\right)  a^{N-1}}_{=\left(  N-2\right)  a^{N}}+2\underbrace{b\left(
N-3\right)  a^{N-2}}_{=-\dfrac{1}{4}\left(  N-3\right)  a^{N}}\right)  x_{0}\\
&  =\dfrac{1}{2^{N-1}}\underbrace{\left(  2\left(  N-1\right)  a^{N-1}+\left(
-\left(  N-2\right)  a^{N-1}\right)  \right)  }_{=Na^{N-1}}x_{1}\\
&  \ \ \ \ \ \ \ \ \ \ -\dfrac{1}{2^{N-1}}\underbrace{\left(  \left(
N-2\right)  a^{N}+2\left(  -\dfrac{1}{4}\left(  N-3\right)  \right)
a^{N}\right)  }_{=\dfrac{1}{2}\left(  N-1\right)  a^{N}}x_{0}\\
&  =\dfrac{1}{2^{N-1}}Na^{N-1}x_{1}-\dfrac{1}{2^{N-1}}\cdot\dfrac{1}{2}\left(
N-1\right)  a^{N}x_{0}=\dfrac{1}{2^{N-1}}\left(  Na^{N-1}x_{1}-\dfrac{1}%
{2}\left(  N-1\right)  a^{N}x_{0}\right) \\
&  =\dfrac{1}{2^{N}}\left(  2Na^{N-1}x_{1}-\left(  N-1\right)  a^{N}%
x_{0}\right)  .
\end{align*}
In other words, (\ref{sol.ps2.2.1.goal}) is proven in Case 3.

Thus we have seen that (\ref{sol.ps2.2.1.goal}) holds in each of the three
Cases 1, 2 and 3. Since these three cases cover all possibilities, this
finishes the proof of (\ref{sol.ps2.2.1.goal}). Hence, we have finished our
proof of (\ref{sol.ps2.2.1.claim}) by strong induction.]
\end{proof}

\begin{remark}
Proving the identity (\ref{sol.ps2.2.1.claim}) by strong induction is a
completely straightforward task. The main difficulty of the exercise is
finding this identity. Linear algebra (specifically, the theory of the Jordan
normal form) gives a \textquotedblleft conceptual\textquotedblright\ way to
derive it (see, e.g., \cite{Fische01} for a very general treatment), but it
can also be experimentally found by computing $x_{2},x_{3},x_{4},x_{5},x_{6}$
directly (using $a^{2}+4b=0$ to rewrite $b$ as $-\dfrac{a^{2}}{4}$, so that
only the variable $a$ appears in the expressions) and guessing the pattern.
\end{remark}

\subsection{Solution to Exercise \ref{exe.ps2.2.2}}

\begin{proof}
[Solution to Exercise \ref{exe.ps2.2.2}.]We shall only solve part
\textbf{(c)}, since the other two parts are its particular cases (for $N=2$
and for $N=3$, respectively).

\textbf{(c)} We define a new sequence $\left(  c_{0},c_{1},c_{2}%
,\ldots\right)  $ recursively by
\begin{align*}
c_{0}  &  =2,\\
c_{1}  &  =a,\ \ \ \ \ \ \ \ \ \ \text{and}\\
c_{n}  &  =ac_{n-1}+bc_{n-2}\ \ \ \ \ \ \ \ \ \ \text{for all }n\geq2.
\end{align*}
(So this sequence $\left(  c_{0},c_{1},c_{2},\ldots\right)  $ is $\left(
a,b\right)  $-recurrent. Its first values are $c_{0}=2$, $c_{1}=a$,
$c_{2}=a^{2}+2b$, $c_{3}=a\left(  a^{2}+3b\right)  $ and $c_{4}=a^{4}%
+4a^{2}b+2b^{2}$. We notice that this sequence depends only on $a$ and $b$.)

We now claim that every $N\in\mathbb{N}$ and $m\in\mathbb{N}$ satisfy%
\begin{equation}
x_{m+2N}=c_{N}x_{m+N}+\left(  -1\right)  ^{N-1}b^{N}x_{m}.
\label{sol.ps2.2.2.c.claim}%
\end{equation}
Once this is proven, we will be done: In fact, (\ref{sol.ps2.2.2.c.claim})
shows that, for every nonnegative integers $N$ and $K$, the sequence $\left(
x_{K},x_{N+K},x_{2N+K},x_{3N+K},\ldots\right)  $ is $\left(  c_{N},\left(
-1\right)  ^{N-1}b^{N}\right)  $-recurrent\footnote{\textit{Proof.} Assume
that we have already proven (\ref{sol.ps2.2.2.c.claim}). Now, for every
nonnegative integers $N$ and $K$, for every $u\geq2$, we have%
\begin{align*}
x_{uN+K}  &  =x_{\left(  u-2\right)  N+K+2N}\ \ \ \ \ \ \ \ \ \ \left(
\text{since }uN+K=\left(  u-2\right)  N+K+2N\right) \\
&  =c_{N}\underbrace{x_{\left(  u-2\right)  N+K+N}}_{=x_{\left(  u-1\right)
N+K}}+\left(  -1\right)  ^{N-1}b^{N}x_{\left(  u-2\right)  N+K}%
\ \ \ \ \ \ \ \ \ \ \left(  \text{by (\ref{sol.ps2.2.2.c.claim}), applied to
}m=\left(  u-2\right)  N+K\right) \\
&  =c_{N}x_{\left(  u-1\right)  N+K}+\left(  -1\right)  ^{N-1}b^{N}x_{\left(
u-2\right)  N+K}.
\end{align*}
In other words, for every nonnegative integers $N$ and $K$, the sequence
$\left(  x_{K},x_{N+K},x_{2N+K},x_{3N+K},\ldots\right)  $ is $\left(
c_{N},\left(  -1\right)  ^{N-1}b^{N}\right)  $-recurrent. Qed.}. Thus, in
order to solve Exercise \ref{exe.ps2.2.2} \textbf{(c)}, we only need to prove
(\ref{sol.ps2.2.2.c.claim}).

[\textit{Proof of (\ref{sol.ps2.2.2.c.claim}):} We shall prove
(\ref{sol.ps2.2.2.c.claim}) by strong induction over $N$. Thus, we fix some
$n\in\mathbb{N}$, and we assume (as our induction hypothesis) that
(\ref{sol.ps2.2.2.c.claim}) holds for every $N<n$ (and, of course, every
$m\in\mathbb{N}$). We need to prove that (\ref{sol.ps2.2.2.c.claim}) holds for
$N=n$ (and every $m\in\mathbb{N}$). In other words, we need to prove that%
\begin{equation}
x_{m+2n}=c_{n}x_{m+n}+\left(  -1\right)  ^{n-1}b^{n}x_{m}
\label{sol.ps2.2.2.c.goal}%
\end{equation}
for every $m\in\mathbb{N}$.

We must be in one of the following three cases:

\textit{Case 1:} We have $n=0$.

\textit{Case 2:} We have $n=1$.

\textit{Case 3:} We have $n\geq2$.

Let us first consider Case 1. In this case, we have $n=0$. Now, let
$m\in\mathbb{N}$. Since $n=0$, we have%
\[
c_{n}x_{m+n}+\left(  -1\right)  ^{n-1}b^{n}x_{m}=\underbrace{c_{0}}%
_{=2}\underbrace{x_{m+0}}_{=x_{m}}+\underbrace{\left(  -1\right)  ^{0-1}%
}_{=-1}\underbrace{b^{0}}_{=1}x_{m}=2x_{m}+\left(  -1\right)  x_{m}=x_{m}.
\]
Compared with%
\begin{align*}
x_{m+2n}  &  =x_{m+2\cdot0}\ \ \ \ \ \ \ \ \ \ \left(  \text{since }n=0\right)
\\
&  =x_{m},
\end{align*}
this yields $x_{m+2n}=c_{n}x_{m+n}+\left(  -1\right)  ^{n-1}b^{n}x_{m}$.
Hence, (\ref{sol.ps2.2.2.c.goal}) is proven in Case 1.

Let us next consider Case 2. In this case, we have $n=1$. Now, let
$m\in\mathbb{N}$. Since $n=1$, we have%
\[
c_{n}x_{m+n}+\left(  -1\right)  ^{n-1}b^{n}x_{m}=\underbrace{c_{1}}%
_{=a}x_{m+1}+\underbrace{\left(  -1\right)  ^{1-1}}_{=1}\underbrace{b^{1}%
}_{=b}x_{m}=ax_{m+1}+bx_{m}.
\]
Compared with%
\begin{align*}
x_{m+2n}  &  =x_{m+2\cdot1}\ \ \ \ \ \ \ \ \ \ \left(  \text{since }n=1\right)
\\
&  =x_{m+2}=ax_{\left(  m+2\right)  -1}+bx_{\left(  m+2\right)  -2}\\
&  \ \ \ \ \ \ \ \ \ \ \left(  \text{since the sequence }\left(  x_{0}%
,x_{1},x_{2},\ldots\right)  \text{ is }\left(  a,b\right)  \text{-recurrent}%
\right) \\
&  =ax_{m+1}+bx_{m},
\end{align*}
this yields $x_{m+2n}=c_{n}x_{m+n}+\left(  -1\right)  ^{n-1}b^{n}x_{m}$.
Hence, (\ref{sol.ps2.2.2.c.goal}) is proven in Case 2.

Let us finally consider Case 3. In this case, we have $n\geq2$. Hence, both
$n-1$ and $n-2$ are nonnegative integers. Moreover, $n-1<n$, so that
(\ref{sol.ps2.2.2.c.claim}) is already proven for $N=n-1$ (by our induction
hypothesis). In other words, we have%
\begin{equation}
x_{m+2\left(  n-1\right)  }=c_{n-1}x_{m+\left(  n-1\right)  }+\left(
-1\right)  ^{\left(  n-1\right)  -1}b^{n-1}x_{m} \label{sol.ps2.2.2.c.hyp1}%
\end{equation}
for every $m\in\mathbb{N}$.

Also, $n-2$ is a nonnegative integer such that $n-2<n$. Hence,
(\ref{sol.ps2.2.1.claim}) is already proven for $N=n-2$ (by our induction
hypothesis). In other words, we have%
\begin{equation}
x_{m+2\left(  n-2\right)  }=c_{n-2}x_{m+\left(  n-2\right)  }+\left(
-1\right)  ^{\left(  n-2\right)  -1}b^{n-2}x_{m} \label{sol.ps2.2.2.c.hyp2}%
\end{equation}
for every $m\in\mathbb{N}$.

Now, fix $m\in\mathbb{N}$. We want to prove (\ref{sol.ps2.2.2.c.goal}). The
recursive definition of the sequence $\left(  c_{0},c_{1},c_{2},\ldots\right)
$ yields
\begin{equation}
c_{n}=ac_{n-1}+bc_{n-2}. \label{sol.ps2.2.2.c.1}%
\end{equation}
Since the sequence $\left(  x_{0},x_{1},x_{2},\ldots\right)  $ is $\left(
a,b\right)  $-recurrent, we have%
\[
x_{m+2}=a\underbrace{x_{\left(  m+2\right)  -1}}_{=x_{m+1}}%
+b\underbrace{x_{\left(  m+2\right)  -2}}_{=x_{m}}=ax_{m+1}+bx_{m},
\]
so that%
\begin{equation}
ax_{m+1}-x_{m+2}=-bx_{m}. \label{sol.ps2.2.2.c.3}%
\end{equation}
But since the sequence $\left(  x_{0},x_{1},x_{2},\ldots\right)  $ is $\left(
a,b\right)  $-recurrent, we also have
\begin{align*}
x_{m+2n}  &  =a\underbrace{x_{\left(  m+2n\right)  -1}}_{\substack{=x_{\left(
m+1\right)  +2\left(  n-1\right)  }\\\text{(since }\left(  m+2n\right)
-1=\left(  m+1\right)  +2\left(  n-1\right)  \text{)}}%
}+b\underbrace{x_{\left(  m+2n\right)  -2}}_{\substack{=x_{\left(  m+2\right)
+2\left(  n-2\right)  }\\\text{(since }\left(  m+2n\right)  -2=\left(
m+2\right)  +2\left(  n-2\right)  \text{)}}}\\
&  \ \ \ \ \ \ \ \ \ \ \left(  \text{since }\underbrace{m}_{\geq
0}+2\underbrace{n}_{\geq2}\geq0+2\cdot2=4\geq2\right) \\
&  =a\underbrace{x_{\left(  m+1\right)  +2\left(  n-1\right)  }}%
_{\substack{=c_{n-1}x_{\left(  m+1\right)  +\left(  n-1\right)  }+\left(
-1\right)  ^{\left(  n-1\right)  -1}b^{n-1}x_{m+1}\\\text{(by
(\ref{sol.ps2.2.2.c.hyp1}), applied to }m+1\text{ instead of }m\text{)}%
}}+b\underbrace{x_{\left(  m+2\right)  +2\left(  n-2\right)  }}%
_{\substack{=c_{n-2}x_{\left(  m+2\right)  +\left(  n-2\right)  }+\left(
-1\right)  ^{\left(  n-2\right)  -1}b^{n-2}x_{m+2}\\\text{(by
(\ref{sol.ps2.2.2.c.hyp2}), applied to }m+2\text{ instead of }m\text{)}}}\\
&  =a\left(  c_{n-1}\underbrace{x_{\left(  m+1\right)  +\left(  n-1\right)  }%
}_{=x_{m+n}}+\underbrace{\left(  -1\right)  ^{\left(  n-1\right)  -1}%
}_{=\left(  -1\right)  ^{n-2}=\left(  -1\right)  ^{n}}b^{n-1}x_{m+1}\right) \\
&  \ \ \ \ \ \ \ \ \ \ +b\left(  c_{n-2}\underbrace{x_{\left(  m+2\right)
+\left(  n-2\right)  }}_{=x_{m+n}}+\underbrace{\left(  -1\right)  ^{\left(
n-2\right)  -1}}_{=\left(  -1\right)  ^{n-3}=\left(  -1\right)  ^{n-1}%
=-\left(  -1\right)  ^{n}}b^{n-2}x_{m+2}\right) \\
&  =a\left(  c_{n-1}x_{m+n}+\left(  -1\right)  ^{n}b^{n-1}x_{m+1}\right)
+b\left(  c_{n-2}x_{m+n}+\left(  -\left(  -1\right)  ^{n}\right)
b^{n-2}x_{m+2}\right) \\
&  =ac_{n-1}x_{m+n}+\left(  -1\right)  ^{n}ab^{n-1}x_{m+1}+bc_{n-2}%
x_{m+n}-\left(  -1\right)  ^{n}\underbrace{bb^{n-2}}_{=b^{n-1}}x_{m+2}\\
&  =ac_{n-1}x_{m+n}+\left(  -1\right)  ^{n}ab^{n-1}x_{m+1}+bc_{n-2}%
x_{m+n}-\left(  -1\right)  ^{n}b^{n-1}x_{m+2}\\
&  =\underbrace{\left(  ac_{n-1}x_{m+n}+bc_{n-2}x_{m+n}\right)  }_{=\left(
ac_{n-1}+bc_{n-2}\right)  x_{m+n}}+\underbrace{\left(  -1\right)  ^{n}%
ab^{n-1}x_{m+1}-\left(  -1\right)  ^{n}b^{n-1}x_{m+2}}_{=\left(  -1\right)
^{n}b^{n-1}\left(  ax_{m+1}-x_{m+2}\right)  }\\
&  =\underbrace{\left(  ac_{n-1}+bc_{n-2}\right)  }_{\substack{=c_{n}%
\\\text{(by (\ref{sol.ps2.2.2.c.1}))}}}x_{m+n}+\left(  -1\right)  ^{n}%
b^{n-1}\underbrace{\left(  ax_{m+1}-x_{m+2}\right)  }_{\substack{=-bx_{m}%
\\\text{(by (\ref{sol.ps2.2.2.c.3}))}}}\\
&  =c_{n}x_{m+n}+\underbrace{\left(  -1\right)  ^{n}b^{n-1}\left(
-bx_{m}\right)  }_{=-\left(  -1\right)  ^{n}b^{n-1}bx_{m}}=c_{n}%
x_{m+n}+\left(  \underbrace{-\left(  -1\right)  ^{n}}_{=\left(  -1\right)
^{n-1}}\underbrace{b^{n-1}b}_{=b^{n}}x_{m}\right) \\
&  =c_{n}x_{m+n}+\left(  -1\right)  ^{n-1}b^{n}x_{m}.
\end{align*}

In other words, (\ref{sol.ps2.2.2.c.goal}) is proven in Case 3.

Thus we have seen that (\ref{sol.ps2.2.2.c.goal}) holds in each of the three
Cases 1, 2 and 3. Since these three cases cover all possibilities, this
finishes the proof of (\ref{sol.ps2.2.2.c.goal}). This finishes our
(inductive) proof of (\ref{sol.ps2.2.2.c.claim}).]

As we know, this solves Exercise \ref{exe.ps2.2.2}.
\end{proof}

\begin{remark}
How on earth could one have come up with my definition of the sequence
$\left(  c_{0},c_{1},c_{2},\ldots\right)  $ in the solution above? One way is
to solve parts \textbf{(a)} and \textbf{(b)} of the exercise first (which can
be solved by applying the equation $x_{n}=ax_{n-1}+bx_{n-2}$ several times),
and then guess that the answer to \textbf{(c)} is a pair of the form $\left(
c,d\right)  =\left(  c_{N},\left(  -1\right)  ^{N-1}b^{N}\right)  $ for some
sequence $\left(  c_{0},c_{1},c_{2},\ldots\right)  $. What remains is finding
this sequence. I believe its entry $c_{3}=a\left(  a^{2}+3b\right)  $ is
particularly telltale.
\end{remark}

\subsection{Solution to Exercise \ref{exe.ps2.2.3}}

\begin{proof}
[Solution to Exercise \ref{exe.ps2.2.3}.]A set $I$ of integers is said to be
\textit{lacunar} if no two elements of $I$ are consecutive (i.e., there exists
no $i\in\mathbb{Z}$ such that both $i$ and $i+1$ belong to $I$). Then,
Exercise \ref{exe.ps2.2.3} asks us to prove that, for every positive integer
$n$,%
\begin{equation}
\text{the number }f_{n}\text{ is the number of lacunar subsets of }\left\{
1,2,\ldots,n-2\right\}  . \label{sol.ps2.2.3.goal}%
\end{equation}

We shall prove (\ref{sol.ps2.2.3.goal}) by strong induction over $n$. Thus, we
let $N$ be a positive integer, and we assume (as the induction hypothesis)
that (\ref{sol.ps2.2.3.goal}) is proven for every $n<N$. We need to prove
(\ref{sol.ps2.2.3.goal}) for $n=N$. In other words, we need to prove that
\begin{equation}
\text{the number }f_{N}\text{ is the number of lacunar subsets of }\left\{
1,2,\ldots,N-2\right\}  . \label{sol.ps2.2.3.goal2}%
\end{equation}

Recall that $N$ is a positive integer. Hence, we are in one of the following
three cases:

\textit{Case 1:} We have $N=1$.

\textit{Case 2:} We have $N=2$.

\textit{Case 3:} We have $N\geq3$.

Let us first consider Case 1. In this case, we have $N=1$. Thus, $f_{N}%
=f_{1}=1$. On the other hand, the number of lacunar subsets of $\left\{
1,2,\ldots,N-2\right\}  $ is $1$ (since the set $\left\{  1,2,\ldots
,\underbrace{N}_{=1}-2\right\}  =\left\{  1,2,\ldots,1-2\right\}
=\varnothing$ has only one subset, and this subset is lacunar). Thus, $f_{N}$
is the number of lacunar subsets of $\left\{  1,2,\ldots,N-2\right\}  $.
Hence, (\ref{sol.ps2.2.3.goal2}) is proven in Case 1.

Case 2 can be dealt with similarly (in this case, the set $\left\{
1,2,\ldots,N-2\right\}  $ is still empty, and $f_{N}$ is still $1$), and is
left to the reader.

We now consider Case 3. In this case, we have $N\geq3$. Hence, $N-1$ and $N-2$
are positive integers. Since $N-1$ is a positive integer and $<N$, we know
that (\ref{sol.ps2.2.3.goal2}) is proven for $n=N-1$ (due to our induction
hypothesis). In other words, the number $f_{N-1}$ is the number of lacunar
subsets of $\left\{  1,2,\ldots,\left(  N-1\right)  -2\right\}  $. In other
words, $f_{N-1}$ is the number of lacunar subsets of $\left\{  1,2,\ldots
,N-3\right\}  $.

Also, since $N-2$ is a positive integer and $<N$, we know that
(\ref{sol.ps2.2.3.goal2}) is proven for $n=N-2$ (due to our induction
hypothesis). In other words, the number $f_{N-2}$ is the number of lacunar
subsets of $\left\{  1,2,\ldots,\left(  N-2\right)  -2\right\}  $. In other
words, $f_{N-2}$ is the number of lacunar subsets of $\left\{  1,2,\ldots
,N-4\right\}  $.

Now, how do the lacunar subsets of $\left\{  1,2,\ldots,N-2\right\}  $ look
like? We say that a lacunar subset of $\left\{  1,2,\ldots,N-2\right\}  $ has
\textit{type 1} if it contains $N-2$, and has \textit{type 2} if it does not.
Now, let us count the lacunar subsets having type 1 and those having type 2:

\begin{itemize}
\item If $I$ is a lacunar subset of $\left\{  1,2,\ldots,N-2\right\}  $ which
has type 1, then it contains $N-2$, and thus cannot contain $N-3$ (because it
is lacunar, i.e., contains no two consecutive integers, but $N-3$ and $N-2$
are two consecutive integers); moreover, $I\setminus\left\{  N-2\right\}  $ is
a lacunar subset of $\left\{  1,2,\ldots,N-4\right\}  $%
\ \ \ \ \footnote{Indeed, it is lacunar because $I$ is lacunar; and it is a
subset of $\left\{  1,2,\ldots,N-4\right\}  $ because $I$ cannot contain
$N-3$.}. Thus, to every lacunar subset $I$ of $\left\{  1,2,\ldots
,N-2\right\}  $ which has type 1, we have assigned a lacunar subset of
$\left\{  1,2,\ldots,N-4\right\}  $ (namely, $I\setminus\left\{  N-2\right\}
$). It is easy to see that this assignment is injective (indeed, if $I$ and
$J$ are two lacunar subsets of $\left\{  1,2,\ldots,N-2\right\}  $ which have
type 1, and if $I\setminus\left\{  N-2\right\}  =J\setminus\left\{
N-2\right\}  $, then $I=J$) and surjective (because whenever $K$ is a lacunar
subset of $\left\{  1,2,\ldots,N-4\right\}  $, the set $K\cup\left\{
N-2\right\}  $ is a lacunar subset of $\left\{  1,2,\ldots,N-2\right\}  $
which has type 1, and this set $K\cup\left\{  N-2\right\}  $ is sent back to
$K$ by our assignment); thus, it is bijective. Hence, we have found a
bijection between the lacunar subsets of $\left\{  1,2,\ldots,N-2\right\}  $
which have type 1 and the lacunar subsets of $\left\{  1,2,\ldots,N-4\right\}
$. Therefore, the number of lacunar subsets of $\left\{  1,2,\ldots
,N-2\right\}  $ which have type 1 equals the number of lacunar subsets of
$\left\{  1,2,\ldots,N-4\right\}  $. But the latter number is $f_{N-2}$.
Hence, the number of lacunar subsets of $\left\{  1,2,\ldots,N-2\right\}  $
which have type 1 is $f_{N-2}$.

\item The lacunar subsets of $\left\{  1,2,\ldots,N-2\right\}  $ which have
type 2 are precisely the lacunar subsets of $\left\{  1,2,\ldots,N-2\right\}
$ which do not contain $N-2$; in other words, they are precisely the lacunar
subsets of $\left\{  1,2,\ldots,N-3\right\}  $. Hence, the number of lacunar
subsets of $\left\{  1,2,\ldots,N-2\right\}  $ which have type 2 is $f_{N-1}$
(since we know that the number of lacunar subsets of $\left\{  1,2,\ldots
,N-3\right\}  $ is $f_{N-1}$).
\end{itemize}

Now, let us summarize. Each of the lacunar subsets of $\left\{  1,2,\ldots
,N-2\right\}  $ either has type 1 or has type 2 (but not both). Hence, the
number of lacunar subsets of $\left\{  1,2,\ldots,N-2\right\}  $ equals%
\begin{align*}
&  \underbrace{\left(  \text{the number of lacunar subsets of }\left\{
1,2,\ldots,N-2\right\}  \text{ which have type 1}\right)  }_{=f_{N-2}}\\
&  \ \ \ \ \ \ \ \ \ \ +\underbrace{\left(  \text{the number of lacunar
subsets of }\left\{  1,2,\ldots,N-2\right\}  \text{ which have type 2}\right)
}_{=f_{N-1}}\\
&  =f_{N-2}+f_{N-1}=f_{N-1}+f_{N-2}\\
&  =f_{N}\ \ \ \ \ \ \ \ \ \ \left(  \text{by the recursion of the Fibonacci
numbers}\right)  .
\end{align*}
This proves (\ref{sol.ps2.2.3.goal2}) in Case 3.

Now, (\ref{sol.ps2.2.3.goal2}) is proven in each of the three Cases 1, 2 and
3. Hence, (\ref{sol.ps2.2.3.goal2}) always holds. This completes our
(inductive) proof of (\ref{sol.ps2.2.3.goal}). Exercise \ref{exe.ps2.2.3} is solved.
\end{proof}

\subsection{Solution to Exercise \ref{exe.ps2.2.S}}

\subsubsection{The solution}

\begin{vershort}
\begin{proof}
[Solution to Exercise \ref{exe.ps2.2.S}.]Let us first explain why the right
hand side of (\ref{eq.exe.2.S}) is well-defined. In fact, this is not obvious,
because if $r=0$, then $r^{n-1-2i}$ might not always make sense (because
$n-1-2i$ can be negative). However, it turns out that $\dbinom{n-1-i}{i}=0$
for every $i\in\left\{  0,1,\ldots,n-1\right\}  $ satisfying $n-1-2i<0$%
\ \ \ \ \footnote{\textit{Proof.} Let $i\in\left\{  0,1,\ldots,n-1\right\}  $
be such that $n-1-2i<0$. Then, $i\leq n-1$ (since $i\in\left\{  0,1,\ldots
,n-1\right\}  $), so that $n-1-i\geq0$. In other words, $n-1-i\in\mathbb{N}$.
But $n-1-i<i$ (since $\left(  n-1-i\right)  -i=n-1-2i<0$). Hence, Proposition
\ref{prop.binom.0} (applied to $n-1-i$ and $i$ instead of $m$ and $n$) yields
$\dbinom{n-1-i}{i}=0$. Qed.}. Hence, we interpret $\dbinom{n-1-i}{i}%
r^{n-1-2i}$ as $0$ whenever $n-1-2i<0$ (even if the term $r^{n-1-2i}$ by
itself is not well-defined), following our convention that any expression of
the form $a\cdot b$ with $a=0$ has to be interpreted as $0$. Thus, the right
hand side of (\ref{eq.exe.2.S}) is well-defined.

Now, we shall prove (\ref{eq.exe.2.S}) by strong induction over $n$. Thus, we
let $N\in\mathbb{N}$, and we assume (as the induction hypothesis) that
(\ref{eq.exe.2.S}) is proven for every $n<N$. We need to prove
(\ref{eq.exe.2.S}) for $n=N$. In other words, we need to prove that
\begin{equation}
c_{N}=\sum_{i=0}^{N-1}\left(  -1\right)  ^{i}\dbinom{N-1-i}{i}r^{N-1-2i}.
\label{sol.ps2.2.S.short.goal}%
\end{equation}

Recall that $N\in\mathbb{N}$. Hence, we are in one of the following three cases:

\textit{Case 1:} We have $N=0$.

\textit{Case 2:} We have $N=1$.

\textit{Case 3:} We have $N\geq2$.

Let us first consider Case 1. In this case, we have $N=0$. Hence, $c_{N}%
=c_{0}=0$. But also, the sum $\sum_{i=0}^{N-1}\left(  -1\right)  ^{i}%
\dbinom{N-1-i}{i}r^{N-1-2i}$ is an empty sum (since $N=0$) and thus equals
$0$. Therefore, both sides of the equality (\ref{sol.ps2.2.S.short.goal})
equal $0$. Hence, the equality (\ref{sol.ps2.2.S.short.goal}) holds. We thus
have proven (\ref{sol.ps2.2.S.short.goal}) in Case 1.

Let us now consider Case 2. In this case, we have $N=1$. Thus, $c_{N}=c_{1}%
=1$. On the other hand, from $N=1$, we obtain%
\begin{align*}
\sum_{i=0}^{N-1}\left(  -1\right)  ^{i}\dbinom{N-1-i}{i}r^{N-1-2i}  &
=\sum_{i=0}^{1-1}\left(  -1\right)  ^{i}\dbinom{1-1-i}{i}r^{1-1-2i}\\
&  =\underbrace{\left(  -1\right)  ^{0}}_{=1}\underbrace{\dbinom{1-1-0}{0}%
}_{=1}\underbrace{r^{1-1-2\cdot0}}_{=r^{0}=1}=1=c_{N}.
\end{align*}
Hence, (\ref{sol.ps2.2.S.short.goal}) is proven in Case 2.

Let us now consider Case 3. In this case, we have $N\geq2$. Therefore, both
$N-1$ and $N-2$ belong to $\mathbb{N}$.

So we know that $N-1$ is an element of $\mathbb{N}$ satisfying $N-1<N$. Hence,
(\ref{eq.exe.2.S}) is proven for $n=N-1$ (by our induction hypothesis). In
other words,%
\begin{equation}
c_{N-1}=\sum_{i=0}^{N-2}\left(  -1\right)  ^{i}\dbinom{N-2-i}{i}r^{N-2-2i}.
\label{sol.ps2.2.S.short.hyp1}%
\end{equation}

Also, $N-2$ is an element of $\mathbb{N}$ satisfying $N-2<N$. Hence,
(\ref{eq.exe.2.S}) is proven for $n=N-2$ (by our induction hypothesis). In
other words,%
\begin{equation}
c_{N-2}=\sum_{i=0}^{N-3}\left(  -1\right)  ^{i}\dbinom{N-3-i}{i}r^{N-3-2i}.
\label{sol.ps2.2.S.short.hyp2}%
\end{equation}

Let us further recall that every positive integer $i$ and every $a\in
\mathbb{Z}$ satisfy%
\begin{equation}
\dbinom{a-1}{i}+\dbinom{a-1}{i-1}=\dbinom{a}{i} \label{sol.ps2.2.S.short.rec1}%
\end{equation}
\footnote{\textit{Proof of (\ref{sol.ps2.2.S.short.rec1}):} Let $i$ be a
positive integer. Let $a\in\mathbb{Z}$. Then, (\ref{eq.binom.rec.m}) (applied
to $m=a$ and $n=i$) shows that $\dbinom{a}{i}=\dbinom{a-1}{i-1}+\dbinom
{a-1}{i}=\dbinom{a-1}{i}+\dbinom{a-1}{i-1}$. This proves
(\ref{sol.ps2.2.S.short.rec1}).}.

Now, $N\geq2$, so that the recursive definition of the sequence $\left(
c_{0},c_{1},c_{2},\ldots\right)  $ yields%
\begin{align*}
c_{N}  &  =r\underbrace{c_{N-1}}_{\substack{=\sum_{i=0}^{N-2}\left(
-1\right)  ^{i}\dbinom{N-2-i}{i}r^{N-2-2i}\\\text{(by
(\ref{sol.ps2.2.S.short.hyp1}))}}}-\underbrace{c_{N-2}}_{\substack{=\sum
_{i=0}^{N-3}\left(  -1\right)  ^{i}\dbinom{N-3-i}{i}r^{N-3-2i}\\\text{(by
(\ref{sol.ps2.2.S.short.hyp2}))}}}\\
&  =r\left(  \sum_{i=0}^{N-2}\left(  -1\right)  ^{i}\dbinom{N-2-i}%
{i}r^{N-2-2i}\right)  -\sum_{i=0}^{N-3}\left(  -1\right)  ^{i}\dbinom
{N-3-i}{i}r^{N-3-2i}\\
&  =\sum_{i=0}^{N-2}\left(  -1\right)  ^{i}\dbinom{N-2-i}{i}%
\underbrace{rr^{N-2-2i}}_{=r^{N-1-2i}}-\underbrace{\sum_{i=0}^{N-3}\left(
-1\right)  ^{i}\dbinom{N-3-i}{i}r^{N-3-2i}}_{\substack{=\sum_{i=1}%
^{N-2}\left(  -1\right)  ^{i-1}\dbinom{N-3-\left(  i-1\right)  }%
{i-1}r^{N-3-2\left(  i-1\right)  }\\\text{(here, we substituted }i-1\text{ for
}i\text{ in the sum)}}}\\
&  =\underbrace{\sum_{i=0}^{N-2}\left(  -1\right)  ^{i}\dbinom{N-2-i}%
{i}r^{N-1-2i}}_{=\left(  -1\right)  ^{0}\dbinom{N-2-0}{0}r^{N-1-2\cdot0}%
+\sum_{i=1}^{N-2}\left(  -1\right)  ^{i}\dbinom{N-2-i}{i}r^{N-1-2i}}\\
&  \ \ \ \ \ \ \ \ \ \ -\sum_{i=1}^{N-2}\underbrace{\left(  -1\right)  ^{i-1}%
}_{=-\left(  -1\right)  ^{i}}\underbrace{\dbinom{N-3-\left(  i-1\right)
}{i-1}}_{=\dbinom{N-2-i}{i-1}}\underbrace{r^{N-3-2\left(  i-1\right)  }%
}_{=r^{N-1-2i}}%
\end{align*}%
\begin{align*}
&  =\left(  -1\right)  ^{0}\underbrace{\dbinom{N-2-0}{0}}_{=1=\dbinom
{N-1-0}{0}}r^{N-1-2\cdot0}\\
&  \ \ \ \ \ \ \ \ \ \ +\underbrace{\sum_{i=1}^{N-2}\left(  -1\right)
^{i}\dbinom{N-2-i}{i}r^{N-1-2i}-\sum_{i=1}^{N-2}\left(  -\left(  -1\right)
^{i}\right)  \dbinom{N-2-i}{i-1}r^{N-1-2i}}_{=\sum_{i=1}^{N-2}\left(
-1\right)  ^{i}\left(  \dbinom{N-2-i}{i}+\dbinom{N-2-i}{i-1}\right)
r^{N-1-2i}}\\
&  =\left(  -1\right)  ^{0}\dbinom{N-1-0}{0}r^{N-1-2\cdot0}\\
&  \ \ \ \ \ \ \ \ \ \ +\sum_{i=1}^{N-2}\left(  -1\right)  ^{i}%
\underbrace{\left(  \dbinom{N-2-i}{i}+\dbinom{N-2-i}{i-1}\right)
}_{\substack{=\dbinom{\left(  N-1-i\right)  -1}{i}+\dbinom{\left(
N-1-i\right)  -1}{i-1}\\=\dbinom{N-1-i}{i}\\\text{(by
(\ref{sol.ps2.2.S.short.rec1}), applied to }a=N-1-i\text{)}}}r^{N-1-2i}\\
&  =\left(  -1\right)  ^{0}\dbinom{N-1-0}{0}r^{N-1-2\cdot0}+\sum_{i=1}%
^{N-2}\left(  -1\right)  ^{i}\dbinom{N-1-i}{i}r^{N-1-2i}\\
&  =\sum_{i=0}^{N-2}\left(  -1\right)  ^{i}\dbinom{N-1-i}{i}r^{N-1-2i}.
\end{align*}
Hence, (\ref{sol.ps2.2.S.short.goal}) is proven in Case 3. We thus have proven
(\ref{sol.ps2.2.S.short.goal}) in all three Cases 1, 2 and 3, so that we
conclude that (\ref{sol.ps2.2.S.short.goal}) always holds. This completes our
proof of (\ref{eq.exe.2.S}).
\end{proof}
\end{vershort}

\begin{verlong}
\begin{proof}
[Solution to Exercise \ref{exe.ps2.2.S}.]Let us first explain why the right
hand side of (\ref{eq.exe.2.S}) is well-defined. In fact, this is not obvious,
because if $r=0$, then $r^{n-1-2i}$ might not always make sense (because
$n-1-2i$ can be negative). However, it turns out that $\dbinom{n-1-i}{i}=0$
for every $i\in\left\{  0,1,\ldots,n-1\right\}  $ satisfying $n-1-2i<0$%
\ \ \ \ \footnote{\textit{Proof.} Let $i\in\left\{  0,1,\ldots,n-1\right\}  $
be such that $n-1-2i<0$. Then, $i\leq n-1$ (since $i\in\left\{  0,1,\ldots
,n-1\right\}  $), so that $n-1-i\geq0$. In other words, $n-1-i\in\mathbb{N}$.
But $n-1-i<i$ (since $\left(  n-1-i\right)  -i=n-1-2i<0$). Hence, Proposition
\ref{prop.binom.0} (applied to $n-1-i$ and $i$ instead of $m$ and $n$) yields
$\dbinom{n-1-i}{i}=0$. Qed.}. Hence, we interpret $\dbinom{n-1-i}{i}%
r^{n-1-2i}$ as $0$ whenever $n-1-2i<0$ (even if the term $r^{n-1-2i}$ by
itself is not well-defined), following our convention that any expression of
the form $a\cdot b$ has to be interpreted as $0$ when $a=0$. Thus, the right
hand side of (\ref{eq.exe.2.S}) is well-defined.

Now, we shall prove (\ref{eq.exe.2.S}) by strong induction over $n$. Thus, we
let $N\in\mathbb{N}$, and we assume (as the induction hypothesis) that
(\ref{eq.exe.2.S}) is proven for every $n<N$. We need to prove
(\ref{eq.exe.2.S}) for $n=N$. In other words, we need to prove that
\begin{equation}
c_{N}=\sum_{i=0}^{N-1}\left(  -1\right)  ^{i}\dbinom{N-1-i}{i}r^{N-1-2i}.
\label{sol.ps2.2.S.goal}%
\end{equation}

Recall that $N\in\mathbb{N}$. Hence, we are in one of the following three cases:

\textit{Case 1:} We have $N=0$.

\textit{Case 2:} We have $N=1$.

\textit{Case 3:} We have $N\geq2$.

Let us first consider Case 1. In this case, we have $N=0$. Hence, $c_{N}%
=c_{0}=0$. But also, the sum $\sum_{i=0}^{N-1}\left(  -1\right)  ^{i}%
\dbinom{N-1-i}{i}r^{N-1-2i}$ is an empty sum (since $N=0$) and thus equals
$0$. Therefore, both sides of the equality (\ref{sol.ps2.2.S.goal}) equal $0$.
Hence, the equality (\ref{sol.ps2.2.S.goal}) holds. We thus have proven
(\ref{sol.ps2.2.S.goal}) in Case 1.

Let us now consider Case 2. In this case, we have $N=1$. Thus, $c_{N}=c_{1}%
=1$. On the other hand, from $N=1$, we obtain%
\begin{align*}
\sum_{i=0}^{N-1}\left(  -1\right)  ^{i}\dbinom{N-1-i}{i}r^{N-1-2i}  &
=\sum_{i=0}^{1-1}\left(  -1\right)  ^{i}\dbinom{1-1-i}{i}r^{1-1-2i}\\
&  =\sum_{i=0}^{0}\left(  -1\right)  ^{i}\dbinom{1-1-i}{i}r^{1-1-2i}\\
&  =\underbrace{\left(  -1\right)  ^{0}}_{=1}\underbrace{\dbinom{1-1-0}{0}%
}_{=1}\underbrace{r^{1-1-2\cdot0}}_{=r^{0}=1}=1=c_{N}.
\end{align*}
Hence, (\ref{sol.ps2.2.S.goal}) is proven in Case 2.

Let us now consider Case 3. In this case, we have $N\geq2$. Therefore, both
$N-1$ and $N-2$ belong to $\mathbb{N}$.

So we know that $N-1$ is an element of $\mathbb{N}$ satisfying $N-1<N$. Hence,
(\ref{eq.exe.2.S}) is proven for $n=N-1$ (by our induction hypothesis). In
other words,%
\[
c_{N-1}=\sum_{i=0}^{\left(  N-1\right)  -1}\left(  -1\right)  ^{i}%
\dbinom{\left(  N-1\right)  -1-i}{i}r^{\left(  N-1\right)  -1-2i}.
\]
This simplifies to%
\begin{equation}
c_{N-1}=\sum_{i=0}^{N-2}\left(  -1\right)  ^{i}\dbinom{N-2-i}{i}r^{N-2-2i}.
\label{sol.ps2.2.S.hyp1}%
\end{equation}

Also, $N-2$ is an element of $\mathbb{N}$ satisfying $N-2<N$. Hence,
(\ref{eq.exe.2.S}) is proven for $n=N-2$ (by our induction hypothesis). In
other words,%
\[
c_{N-2}=\sum_{i=0}^{\left(  N-2\right)  -1}\left(  -1\right)  ^{i}%
\dbinom{\left(  N-2\right)  -1-i}{i}r^{\left(  N-2\right)  -1-2i}.
\]
This simplifies to%
\begin{equation}
c_{N-2}=\sum_{i=0}^{N-3}\left(  -1\right)  ^{i}\dbinom{N-3-i}{i}r^{N-3-2i}.
\label{sol.ps2.2.S.hyp2}%
\end{equation}

Let us further recall that every positive integer $i$ and every $a\in
\mathbb{Z}$ satisfy%
\begin{equation}
\dbinom{a}{i}=\dbinom{a-1}{i}+\dbinom{a-1}{i-1} \label{sol.ps2.2.S.rec1}%
\end{equation}
\footnote{\textit{Proof of (\ref{sol.ps2.2.S.rec1}):} Let $i$ be a positive
integer. Let $a\in\mathbb{Z}$. Then, (\ref{eq.binom.rec.m}) (applied to $m=a$
and $n=i$) shows that $\dbinom{a}{i}=\dbinom{a-1}{i-1}+\dbinom{a-1}{i}%
=\dbinom{a-1}{i}+\dbinom{a-1}{i-1}$. This proves (\ref{sol.ps2.2.S.rec1}).}.

\begin{noncompile}
Here is a proof of (\ref{sol.ps2.2.S.rec1}) that I have written back before I
had written up the section about binomial coefficients:

[\textit{Proof of (\ref{sol.ps2.2.S.rec1}):} This is the famous recurrence
relation of Pascal's triangle, saying that every number in Pascal's triangle
equals the sum of the two numbers above. However, Pascal's triangle is usually
only drawn to contain the binomial coefficients $\dbinom{a}{i}$ with $0\leq
i\leq a$ (or else it would not be a triangle), so you might not be aware that
it holds for all $a\in\mathbb{Z}$. Either way, the proof is simple:%
\begin{align*}
&  \underbrace{\dbinom{a-1}{i}}_{\substack{=\dfrac{\left(  a-1\right)  \left(
a-2\right)  \cdots\left(  a-i\right)  }{i!}=\dfrac{\left(  a-1\right)  \left(
a-2\right)  \cdots\left(  a-i+1\right)  }{\left(  i-1\right)  !}\cdot
\dfrac{a-i}{i}\\\text{(since }i!=\left(  i-1\right)  !\cdot i\text{)}%
}}+\underbrace{\dbinom{a-1}{i-1}}_{=\dfrac{\left(  a-1\right)  \left(
a-2\right)  \cdots\left(  a-i+1\right)  }{\left(  i-1\right)  !}}\\
&  =\dfrac{\left(  a-1\right)  \left(  a-2\right)  \cdots\left(  a-i+1\right)
}{\left(  i-1\right)  !}\cdot\dfrac{a-i}{i}+\dfrac{\left(  a-1\right)  \left(
a-2\right)  \cdots\left(  a-i+1\right)  }{\left(  i-1\right)  !}\\
&  =\dfrac{\left(  a-1\right)  \left(  a-2\right)  \cdots\left(  a-i+1\right)
}{\left(  i-1\right)  !}\underbrace{\left(  \dfrac{a-i}{i}+1\right)
}_{=\dfrac{a}{i}}=\dfrac{\left(  a-1\right)  \left(  a-2\right)  \cdots\left(
a-i+1\right)  }{\left(  i-1\right)  !}\cdot\dfrac{a}{i}\\
&  =\dfrac{a\left(  a-1\right)  \left(  a-2\right)  \cdots\left(
a-i+1\right)  }{i!}\ \ \ \ \ \ \ \ \ \ \left(  \text{since }\left(
i-1\right)  !\cdot i=i!\right) \\
&  =\dbinom{a}{i}.
\end{align*}
This proves (\ref{sol.ps2.2.S.rec1}).]
\end{noncompile}

Now, $N\geq2$, so that the recursive definition of the sequence $\left(
c_{0},c_{1},c_{2},\ldots\right)  $ yields%
\begin{align*}
c_{N}  &  =r\underbrace{c_{N-1}}_{\substack{=\sum_{i=0}^{N-2}\left(
-1\right)  ^{i}\dbinom{N-2-i}{i}r^{N-2-2i}\\\text{(by (\ref{sol.ps2.2.S.hyp1}%
))}}}-\underbrace{c_{N-2}}_{\substack{=\sum_{i=0}^{N-3}\left(  -1\right)
^{i}\dbinom{N-3-i}{i}r^{N-3-2i}\\\text{(by (\ref{sol.ps2.2.S.hyp2}))}}}\\
&  =r\left(  \sum_{i=0}^{N-2}\left(  -1\right)  ^{i}\dbinom{N-2-i}%
{i}r^{N-2-2i}\right)  -\sum_{i=0}^{N-3}\left(  -1\right)  ^{i}\dbinom
{N-3-i}{i}r^{N-3-2i}\\
&  =\sum_{i=0}^{N-2}\left(  -1\right)  ^{i}\dbinom{N-2-i}{i}%
\underbrace{rr^{N-2-2i}}_{\substack{=r^{\left(  N-2-2i\right)  +1}%
\\=r^{N-1-2i}}}-\underbrace{\sum_{i=0}^{N-3}\left(  -1\right)  ^{i}%
\dbinom{N-3-i}{i}r^{N-3-2i}}_{\substack{=\sum_{i=1}^{N-2}\left(  -1\right)
^{i-1}\dbinom{N-3-\left(  i-1\right)  }{i-1}r^{N-3-2\left(  i-1\right)
}\\\text{(here, we substituted }i-1\text{ for }i\text{ in the sum)}}}\\
&  =\underbrace{\sum_{i=0}^{N-2}\left(  -1\right)  ^{i}\dbinom{N-2-i}%
{i}r^{N-1-2i}}_{=\left(  -1\right)  ^{0}\dbinom{N-2-0}{0}r^{N-1-2\cdot0}%
+\sum_{i=1}^{N-2}\left(  -1\right)  ^{i}\dbinom{N-2-i}{i}r^{N-1-2i}}\\
&  \ \ \ \ \ \ \ \ \ \ -\sum_{i=1}^{N-2}\underbrace{\left(  -1\right)  ^{i-1}%
}_{=-\left(  -1\right)  ^{i}}\underbrace{\dbinom{N-3-\left(  i-1\right)
}{i-1}}_{=\dbinom{N-2-i}{i-1}}\underbrace{r^{N-3-2\left(  i-1\right)  }%
}_{=r^{N-1-2i}}%
\end{align*}%
\begin{align*}
&  =\left(  -1\right)  ^{0}\underbrace{\dbinom{N-2-0}{0}}_{=1=\dbinom
{N-1-0}{0}}r^{N-1-2\cdot0}+\sum_{i=1}^{N-2}\left(  -1\right)  ^{i}%
\dbinom{N-2-i}{i}r^{N-1-2i}\\
&  \ \ \ \ \ \ \ \ \ \ -\underbrace{\sum_{i=1}^{N-2}\left(  -\left(
-1\right)  ^{i}\right)  \dbinom{N-2-i}{i-1}r^{N-1-2i}}_{=-\sum_{i=1}%
^{N-2}\left(  -1\right)  ^{i}\dbinom{N-2-i}{i-1}r^{N-1-2i}}\\
&  =\left(  -1\right)  ^{0}\dbinom{N-1-0}{0}r^{N-1-2\cdot0}+\sum_{i=1}%
^{N-2}\left(  -1\right)  ^{i}\dbinom{N-2-i}{i}r^{N-1-2i}\\
&  \ \ \ \ \ \ \ \ \ \ -\left(  -\sum_{i=1}^{N-2}\left(  -1\right)
^{i}\dbinom{N-2-i}{i-1}r^{N-1-2i}\right) \\
&  =\left(  -1\right)  ^{0}\dbinom{N-1-0}{0}r^{N-1-2\cdot0}\\
&  \ \ \ \ \ \ \ \ \ \ +\underbrace{\sum_{i=1}^{N-2}\left(  -1\right)
^{i}\dbinom{N-2-i}{i}r^{N-1-2i}+\sum_{i=1}^{N-2}\left(  -1\right)  ^{i}%
\dbinom{N-2-i}{i-1}r^{N-1-2i}}_{=\sum_{i=1}^{N-2}\left(  -1\right)
^{i}\left(  \dbinom{N-2-i}{i}+\dbinom{N-2-i}{i-1}\right)  r^{N-1-2i}}\\
&  =\left(  -1\right)  ^{0}\dbinom{N-1-0}{0}r^{N-1-2\cdot0}\\
&  \ \ \ \ \ \ \ \ \ \ +\sum_{i=1}^{N-2}\left(  -1\right)  ^{i}\left(
\dbinom{N-2-i}{i}+\dbinom{N-2-i}{i-1}\right)  r^{N-1-2i}.
\end{align*}
Comparing this with%
\begin{align*}
&  \sum_{i=0}^{N-2}\left(  -1\right)  ^{i}\dbinom{N-1-i}{i}r^{N-1-2i}\\
&  =\left(  -1\right)  ^{0}\dbinom{N-1-0}{0}r^{N-1-2\cdot0}\\
&  \ \ \ \ \ \ \ \ \ \ +\sum_{i=1}^{N-2}\left(  -1\right)  ^{i}%
\underbrace{\dbinom{N-1-i}{i}}_{\substack{=\dbinom{\left(  N-1-i\right)
-1}{i}+\dbinom{\left(  N-1-i\right)  -1}{i-1}\\\text{(by
(\ref{sol.ps2.2.S.rec1}), applied to }a=N-1-i\text{)}}}r^{N-1-2i}\\
&  =\left(  -1\right)  ^{0}\dbinom{N-1-0}{0}r^{N-1-2\cdot0}\\
&  \ \ \ \ \ \ \ \ \ \ +\sum_{i=1}^{N-2}\left(  -1\right)  ^{i}\left(
\underbrace{\dbinom{\left(  N-1-i\right)  -1}{i}}_{=\dbinom{N-2-i}{i}%
}+\underbrace{\dbinom{\left(  N-1-i\right)  -1}{i-1}}_{=\dbinom{N-2-i}{i-1}%
}\right)  r^{N-1-2i}\\
&  =\left(  -1\right)  ^{0}\dbinom{N-1-0}{0}r^{N-1-2\cdot0}\\
&  \ \ \ \ \ \ \ \ \ \ +\sum_{i=1}^{N-2}\left(  -1\right)  ^{i}\left(
\dbinom{N-2-i}{i}+\dbinom{N-2-i}{i-1}\right)  r^{N-1-2i},
\end{align*}
we obtain%
\[
c_{N}=\sum_{i=0}^{N-1}\left(  -1\right)  ^{i}\dbinom{N-1-i}{i}r^{N-1-2i}.
\]
Hence, (\ref{sol.ps2.2.S.goal}) is proven in Case 3. We thus have proven
(\ref{sol.ps2.2.S.goal}) in all three Cases 1, 2 and 3, so that we conclude
that (\ref{sol.ps2.2.S.goal}) always holds. This completes our proof of
(\ref{eq.exe.2.S}).
\end{proof}
\end{verlong}

\subsubsection{A corollary}

Having solved Exercise \ref{exe.ps2.2.S}, let us show an identity for binomial
coefficients that follows from it:

\begin{corollary}
\label{cor.ps2.2.S.-1}For each $n\in\mathbb{N}$, we have%
\[
\sum_{i=0}^{n-1}\left(  -1\right)  ^{i}\dbinom{n-1-i}{i}=\left(  -1\right)
^{n}\cdot%
\begin{cases}
0, & \text{if }n\equiv0\operatorname{mod}3;\\
-1, & \text{if }n\equiv1\operatorname{mod}3;\\
1, & \text{if }n\equiv2\operatorname{mod}3
\end{cases}
.
\]
{}
\end{corollary}

\begin{proof}
[Proof of Corollary \ref{cor.ps2.2.S.-1}.]Define $r\in\mathbb{Z}$ by $r=-1$.
Define a sequence $\left(  c_{0},c_{1},c_{2},\ldots\right)  $ of integers as
in Exercise \ref{exe.ps2.2.S}. (Thus, $c_{0}=0$, $c_{1}=1$ and $c_{n}%
=rc_{n-1}-c_{n-2}$ for all $n\geq2$.)

\begin{statement}
\textit{Observation 1:} We have%
\[
c_{n}=%
\begin{cases}
0, & \text{if }n\equiv0\operatorname{mod}3;\\
1, & \text{if }n\equiv1\operatorname{mod}3;\\
-1, & \text{if }n\equiv2\operatorname{mod}3
\end{cases}
\]
for each $n\in\mathbb{N}$.
\end{statement}

[\textit{Proof of Observation 1:} We shall prove Observation 1 by strong
induction over $n$:

\textit{Induction step:} Let $N\in\mathbb{N}$. Assume that Observation 1 holds
whenever $n<N$. We must now prove that Observation 1 holds for $n=N$. In other
words, we must prove that%
\begin{equation}
c_{N}=%
\begin{cases}
0, & \text{if }N\equiv0\operatorname{mod}3;\\
1, & \text{if }N\equiv1\operatorname{mod}3;\\
-1, & \text{if }N\equiv2\operatorname{mod}3
\end{cases}
. \label{pf.cor.ps2.2.S.-1.o1.pf.goal}%
\end{equation}

\begin{vershort}
If $N<2$, then this can be checked directly (since $c_{0}=0$ and $c_{1}=1$).
Thus, for the rest of this proof, we WLOG assume that we don't have $N<2$.
Hence, $N\geq2$; therefore, $N-1\in\mathbb{N}$ and $N-2\in\mathbb{N}$.
\end{vershort}

\begin{verlong}
If $N<2$, then this is easy to check\footnote{\textit{Proof.} Assume that
$N<2$. We must prove (\ref{pf.cor.ps2.2.S.-1.o1.pf.goal}).
\par
We have $N<2$. Thus, $N\leq1$ (since $N$ is an integer), so that $N\in\left\{
0,1\right\}  $ (since $N\in\mathbb{N}$). In other words, either $N=0$ or
$N=1$. Hence, we are in one of the following two cases:
\par
\textit{Case 1:} We have $N=0$.
\par
\textit{Case 2:} We have $N=1$.
\par
Let us first consider Case 1. In this case, we have $N=0$. Hence, $c_{N}%
=c_{0}=0$. But $N=0\equiv0\operatorname{mod}3$, so that $%
\begin{cases}
0, & \text{if }N\equiv0\operatorname{mod}3;\\
1, & \text{if }N\equiv1\operatorname{mod}3;\\
-1, & \text{if }N\equiv2\operatorname{mod}3
\end{cases}
=0$. Comparing this with $c_{N}=0$, we obtain $c_{N}=%
\begin{cases}
0, & \text{if }N\equiv0\operatorname{mod}3;\\
1, & \text{if }N\equiv1\operatorname{mod}3;\\
-1, & \text{if }N\equiv2\operatorname{mod}3
\end{cases}
$. Thus, (\ref{pf.cor.ps2.2.S.-1.o1.pf.goal}) is proven in Case 1.
\par
Let us now consider Case 2. In this case, we have $N=1$. Hence, $c_{N}%
=c_{1}=1$. But $N=1\equiv1\operatorname{mod}3$, so that $%
\begin{cases}
0, & \text{if }N\equiv0\operatorname{mod}3;\\
1, & \text{if }N\equiv1\operatorname{mod}3;\\
-1, & \text{if }N\equiv2\operatorname{mod}3
\end{cases}
=1$. Comparing this with $c_{N}=1$, we obtain $c_{N}=%
\begin{cases}
0, & \text{if }N\equiv0\operatorname{mod}3;\\
1, & \text{if }N\equiv1\operatorname{mod}3;\\
-1, & \text{if }N\equiv2\operatorname{mod}3
\end{cases}
$. Thus, (\ref{pf.cor.ps2.2.S.-1.o1.pf.goal}) is proven in Case 2.
\par
We have now proven (\ref{pf.cor.ps2.2.S.-1.o1.pf.goal}) in each of the two
Cases 1 and 2. Since these two Cases cover all possibilities, we thus conclude
that (\ref{pf.cor.ps2.2.S.-1.o1.pf.goal}) always holds (under the assumption
that $N<2$). Qed.}. Thus, for the rest of this proof, we can WLOG assume that
we don't have $N<2$. Assume this.

We have $N\geq2$ (since we don't have $N<2$). Thus, $N-1\in\mathbb{N}$ and
$N-2\in\mathbb{N}$.
\end{verlong}

But the recursive definition of the sequence $\left(  c_{0},c_{1},c_{2}%
,\ldots\right)  $ yields $c_{n}=rc_{n-1}-c_{n-2}$ for all $n\geq2$. Applying
this to $n=N$, we find%
\begin{equation}
c_{N}=\underbrace{r}_{=-1}c_{N-1}-c_{N-2}=\left(  -1\right)  c_{N-1}%
-c_{N-2}=-c_{N-1}-c_{N-2}. \label{pf.cor.ps2.2.S.-1.o1.pf.rec}%
\end{equation}

We have assumed that Observation 1 holds whenever $n<N$. Thus, we can apply
Observation 1 to $n=N-1$ (since $N-1<N$ and $N-1\in\mathbb{N}$). We thus
conclude that%
\begin{equation}
c_{N-1}=%
\begin{cases}
0, & \text{if }N-1\equiv0\operatorname{mod}3;\\
1, & \text{if }N-1\equiv1\operatorname{mod}3;\\
-1, & \text{if }N-1\equiv2\operatorname{mod}3
\end{cases}
. \label{pf.cor.ps2.2.S.-1.o1.pf.1}%
\end{equation}

\begin{vershort}
Similarly, we obtain%
\begin{equation}
c_{N-2}=%
\begin{cases}
0, & \text{if }N-2\equiv0\operatorname{mod}3;\\
1, & \text{if }N-2\equiv1\operatorname{mod}3;\\
-1, & \text{if }N-2\equiv2\operatorname{mod}3
\end{cases}
. \label{pf.cor.ps2.2.S.-1.o1.pf.short.2}%
\end{equation}

\end{vershort}

\begin{verlong}
We have assumed that Observation 1 holds whenever $n<N$. Thus, we can apply
Observation 1 to $n=N-2$ (since $N-2<N$ and $N-2\in\mathbb{N}$). We thus
conclude that%
\begin{equation}
c_{N-2}=%
\begin{cases}
0, & \text{if }N-2\equiv0\operatorname{mod}3;\\
1, & \text{if }N-2\equiv1\operatorname{mod}3;\\
-1, & \text{if }N-2\equiv2\operatorname{mod}3
\end{cases}
. \label{pf.cor.ps2.2.S.-1.o1.pf.2}%
\end{equation}

\end{verlong}

The remainder obtained when $N$ is divided by $3$ must be either $0$ or $1$ or
$2$. Hence, we must have either $N\equiv0\operatorname{mod}3$ or
$N\equiv1\operatorname{mod}3$ or $N\equiv2\operatorname{mod}3$. Thus, we are
in one of the following three cases:

\textit{Case 1:} We have $N\equiv0\operatorname{mod}3$.

\textit{Case 2:} We have $N\equiv1\operatorname{mod}3$.

\textit{Case 3:} We have $N\equiv2\operatorname{mod}3$.

\begin{vershort}
We shall only consider Case 1 (and leave the two other cases, which are
analogous, to the reader). In this case, we have $N\equiv0\operatorname{mod}%
3$. Hence, $\underbrace{N}_{\equiv0\operatorname{mod}3}-1\equiv0-1\equiv
2\operatorname{mod}3$. Thus, (\ref{pf.cor.ps2.2.S.-1.o1.pf.1}) simplifies to
$c_{N-1}=-1$. Also, $\underbrace{N}_{\equiv0\operatorname{mod}3}%
-2\equiv0-2\equiv1\operatorname{mod}3$. Thus,
(\ref{pf.cor.ps2.2.S.-1.o1.pf.short.2}) simplifies to $c_{N-2}=1$. Now,
(\ref{pf.cor.ps2.2.S.-1.o1.pf.rec}) yields%
\[
c_{N}=-\underbrace{c_{N-1}}_{=-1}-\underbrace{c_{N-2}}_{=1}=-\left(
-1\right)  -1=0.
\]
Comparing this with%
\[%
\begin{cases}
0, & \text{if }N\equiv0\operatorname{mod}3;\\
1, & \text{if }N\equiv1\operatorname{mod}3;\\
-1, & \text{if }N\equiv2\operatorname{mod}3
\end{cases}
=0\ \ \ \ \ \ \ \ \ \ \left(  \text{since }N\equiv0\operatorname{mod}3\right)
,
\]
we obtain $c_{N}=%
\begin{cases}
0, & \text{if }N\equiv0\operatorname{mod}3;\\
1, & \text{if }N\equiv1\operatorname{mod}3;\\
-1, & \text{if }N\equiv2\operatorname{mod}3
\end{cases}
$. Thus, (\ref{pf.cor.ps2.2.S.-1.o1.pf.goal}) is proven in Case 1. As we have
said, the proof in the other two cases is analogous. Hence,
(\ref{pf.cor.ps2.2.S.-1.o1.pf.goal}) is proven.
\end{vershort}

\begin{verlong}
Let us first consider Case 1. In this case, we have $N\equiv
0\operatorname{mod}3$. Hence, $\underbrace{N}_{\equiv0\operatorname{mod}%
3}-1\equiv0-1\equiv2\operatorname{mod}3$. Thus,
(\ref{pf.cor.ps2.2.S.-1.o1.pf.1}) yields%
\[
c_{N-1}=%
\begin{cases}
0, & \text{if }N-1\equiv0\operatorname{mod}3;\\
1, & \text{if }N-1\equiv1\operatorname{mod}3;\\
-1, & \text{if }N-1\equiv2\operatorname{mod}3
\end{cases}
=-1
\]
(since $N-1\equiv2\operatorname{mod}3$). Also, $\underbrace{N}_{\equiv
0\operatorname{mod}3}-2\equiv0-2\equiv1\operatorname{mod}3$. Thus,
(\ref{pf.cor.ps2.2.S.-1.o1.pf.2}) yields%
\[
c_{N-2}=%
\begin{cases}
0, & \text{if }N-2\equiv0\operatorname{mod}3;\\
1, & \text{if }N-2\equiv1\operatorname{mod}3;\\
-1, & \text{if }N-2\equiv2\operatorname{mod}3
\end{cases}
=1
\]
(since $N-2\equiv1\operatorname{mod}3$). Now,
(\ref{pf.cor.ps2.2.S.-1.o1.pf.rec}) yields%
\[
c_{N}=-\underbrace{c_{N-1}}_{=-1}-\underbrace{c_{N-2}}_{=1}=-\left(
-1\right)  -1=0.
\]
Comparing this with%
\[%
\begin{cases}
0, & \text{if }N\equiv0\operatorname{mod}3;\\
1, & \text{if }N\equiv1\operatorname{mod}3;\\
-1, & \text{if }N\equiv2\operatorname{mod}3
\end{cases}
=0\ \ \ \ \ \ \ \ \ \ \left(  \text{since }N\equiv0\operatorname{mod}3\right)
,
\]
we obtain $c_{N}=%
\begin{cases}
0, & \text{if }N\equiv0\operatorname{mod}3;\\
1, & \text{if }N\equiv1\operatorname{mod}3;\\
-1, & \text{if }N\equiv2\operatorname{mod}3
\end{cases}
$. Thus, (\ref{pf.cor.ps2.2.S.-1.o1.pf.goal}) is proven in Case 1.

Let us next consider Case 2. In this case, we have $N\equiv1\operatorname{mod}%
3$. Hence, $\underbrace{N}_{\equiv1\operatorname{mod}3}-1\equiv1-1\equiv
0\operatorname{mod}3$. Thus, (\ref{pf.cor.ps2.2.S.-1.o1.pf.1}) yields%
\[
c_{N-1}=%
\begin{cases}
0, & \text{if }N-1\equiv0\operatorname{mod}3;\\
1, & \text{if }N-1\equiv1\operatorname{mod}3;\\
-1, & \text{if }N-1\equiv2\operatorname{mod}3
\end{cases}
=0
\]
(since $N-1\equiv0\operatorname{mod}3$). Also, $\underbrace{N}_{\equiv
1\operatorname{mod}3}-2\equiv1-2\equiv2\operatorname{mod}3$. Thus,
(\ref{pf.cor.ps2.2.S.-1.o1.pf.2}) yields%
\[
c_{N-2}=%
\begin{cases}
0, & \text{if }N-2\equiv0\operatorname{mod}3;\\
1, & \text{if }N-2\equiv1\operatorname{mod}3;\\
-1, & \text{if }N-2\equiv2\operatorname{mod}3
\end{cases}
=-1
\]
(since $N-2\equiv2\operatorname{mod}3$). Now,
(\ref{pf.cor.ps2.2.S.-1.o1.pf.rec}) yields%
\[
c_{N}=-\underbrace{c_{N-1}}_{=0}-\underbrace{c_{N-2}}_{=-1}=-0-\left(
-1\right)  =1.
\]
Comparing this with%
\[%
\begin{cases}
0, & \text{if }N\equiv0\operatorname{mod}3;\\
1, & \text{if }N\equiv1\operatorname{mod}3;\\
-1, & \text{if }N\equiv2\operatorname{mod}3
\end{cases}
=1\ \ \ \ \ \ \ \ \ \ \left(  \text{since }N\equiv1\operatorname{mod}3\right)
,
\]
we obtain $c_{N}=%
\begin{cases}
0, & \text{if }N\equiv0\operatorname{mod}3;\\
1, & \text{if }N\equiv1\operatorname{mod}3;\\
-1, & \text{if }N\equiv2\operatorname{mod}3
\end{cases}
$. Thus, (\ref{pf.cor.ps2.2.S.-1.o1.pf.goal}) is proven in Case 2.

Let us finally consider Case 3. In this case, we have $N\equiv
2\operatorname{mod}3$. Hence, $\underbrace{N}_{\equiv2\operatorname{mod}%
3}-1\equiv2-1=1\operatorname{mod}3$. Thus, (\ref{pf.cor.ps2.2.S.-1.o1.pf.1})
yields%
\[
c_{N-1}=%
\begin{cases}
0, & \text{if }N-1\equiv0\operatorname{mod}3;\\
1, & \text{if }N-1\equiv1\operatorname{mod}3;\\
-1, & \text{if }N-1\equiv2\operatorname{mod}3
\end{cases}
=1
\]
(since $N-1\equiv1\operatorname{mod}3$). Also, $\underbrace{N}_{\equiv
2\operatorname{mod}3}-2\equiv2-2=0\operatorname{mod}3$. Thus,
(\ref{pf.cor.ps2.2.S.-1.o1.pf.2}) yields%
\[
c_{N-2}=%
\begin{cases}
0, & \text{if }N-2\equiv0\operatorname{mod}3;\\
1, & \text{if }N-2\equiv1\operatorname{mod}3;\\
-1, & \text{if }N-2\equiv2\operatorname{mod}3
\end{cases}
=0
\]
(since $N-2\equiv0\operatorname{mod}3$). Now,
(\ref{pf.cor.ps2.2.S.-1.o1.pf.rec}) yields%
\[
c_{N}=-\underbrace{c_{N-1}}_{=1}-\underbrace{c_{N-2}}_{=0}=-1-0=-1.
\]
Comparing this with%
\[%
\begin{cases}
0, & \text{if }N\equiv0\operatorname{mod}3;\\
1, & \text{if }N\equiv1\operatorname{mod}3;\\
-1, & \text{if }N\equiv2\operatorname{mod}3
\end{cases}
=-1\ \ \ \ \ \ \ \ \ \ \left(  \text{since }N\equiv2\operatorname{mod}%
3\right)  ,
\]
we obtain $c_{N}=%
\begin{cases}
0, & \text{if }N\equiv0\operatorname{mod}3;\\
1, & \text{if }N\equiv1\operatorname{mod}3;\\
-1, & \text{if }N\equiv2\operatorname{mod}3
\end{cases}
$. Thus, (\ref{pf.cor.ps2.2.S.-1.o1.pf.goal}) is proven in Case 3.

We have thus proven (\ref{pf.cor.ps2.2.S.-1.o1.pf.goal}) in each of the three
Cases 1, 2 and 3. Since these three Cases cover all possibilities, we thus
conclude that (\ref{pf.cor.ps2.2.S.-1.o1.pf.goal}) always holds.
\end{verlong}

In other words, Observation 1 holds for $n=N$. This completes the induction
step. Thus, Observation 1 is proven.]

Now, let $n\in\mathbb{N}$. Clearly, $\left(  -1\right)  ^{n-1}\left(
-1\right)  ^{n-1}=\left(  \underbrace{\left(  -1\right)  \left(  -1\right)
}_{=1}\right)  ^{n-1}=1^{n-1}=1$.

Observation 1 yields%
\[
c_{n}=%
\begin{cases}
0, & \text{if }n\equiv0\operatorname{mod}3;\\
1, & \text{if }n\equiv1\operatorname{mod}3;\\
-1, & \text{if }n\equiv2\operatorname{mod}3
\end{cases}
.
\]
Thus,%
\begin{align}
-c_{n}  &  =-%
\begin{cases}
0, & \text{if }n\equiv0\operatorname{mod}3;\\
1, & \text{if }n\equiv1\operatorname{mod}3;\\
-1, & \text{if }n\equiv2\operatorname{mod}3
\end{cases}
=%
\begin{cases}
-0, & \text{if }n\equiv0\operatorname{mod}3;\\
-1, & \text{if }n\equiv1\operatorname{mod}3;\\
-\left(  -1\right)  , & \text{if }n\equiv2\operatorname{mod}3
\end{cases}
\nonumber\\
&  =%
\begin{cases}
0, & \text{if }n\equiv0\operatorname{mod}3;\\
-1, & \text{if }n\equiv1\operatorname{mod}3;\\
1, & \text{if }n\equiv2\operatorname{mod}3
\end{cases}
. \label{pf.cor.ps2.2.S.-cn}%
\end{align}

Exercise \ref{exe.ps2.2.S} yields%
\begin{align*}
c_{n}  &  =\sum_{i=0}^{n-1}\left(  -1\right)  ^{i}\dbinom{n-1-i}%
{i}\underbrace{r^{n-1-2i}}_{\substack{=\left(  -1\right)  ^{n-1-2i}%
\\\text{(since }r=-1\text{)}}}\\
&  =\sum_{i=0}^{n-1}\left(  -1\right)  ^{i}\dbinom{n-1-i}{i}%
\underbrace{\left(  -1\right)  ^{n-1-2i}}_{\substack{=\left(  -1\right)
^{n-1}\\\text{(since }n-1-2i\equiv n-1\operatorname{mod}2\text{)}}}\\
&  =\sum_{i=0}^{n-1}\left(  -1\right)  ^{i}\dbinom{n-1-i}{i}\left(  -1\right)
^{n-1}.
\end{align*}
Multiplying both sides of this equality by $\left(  -1\right)  ^{n-1}$, we
obtain%
\begin{align*}
\left(  -1\right)  ^{n-1}c_{n}  &  =\left(  -1\right)  ^{n-1}\sum_{i=0}%
^{n-1}\left(  -1\right)  ^{i}\dbinom{n-1-i}{i}\left(  -1\right)  ^{n-1}\\
&  =\sum_{i=0}^{n-1}\left(  -1\right)  ^{i}\dbinom{n-1-i}{i}%
\underbrace{\left(  -1\right)  ^{n-1}\left(  -1\right)  ^{n-1}}_{=1}\\
&  =\sum_{i=0}^{n-1}\left(  -1\right)  ^{i}\dbinom{n-1-i}{i}.
\end{align*}
Therefore,%
\begin{align*}
\sum_{i=0}^{n-1}\left(  -1\right)  ^{i}\dbinom{n-1-i}{i}  &
=\underbrace{\left(  -1\right)  ^{n-1}}_{=-\left(  -1\right)  ^{n}}%
c_{n}=-\left(  -1\right)  ^{n}c_{n}\\
&  =\left(  -1\right)  ^{n}\cdot\underbrace{\left(  -c_{n}\right)
}_{\substack{=%
\begin{cases}
0, & \text{if }n\equiv0\operatorname{mod}3;\\
-1, & \text{if }n\equiv1\operatorname{mod}3;\\
1, & \text{if }n\equiv2\operatorname{mod}3
\end{cases}
\\\text{(by (\ref{pf.cor.ps2.2.S.-cn}))}}}\\
&  =\left(  -1\right)  ^{n}\cdot%
\begin{cases}
0, & \text{if }n\equiv0\operatorname{mod}3;\\
-1, & \text{if }n\equiv1\operatorname{mod}3;\\
1, & \text{if }n\equiv2\operatorname{mod}3
\end{cases}
.
\end{align*}
This proves Corollary \ref{cor.ps2.2.S.-1}.
\end{proof}

Let us remark that Corollary \ref{cor.ps2.2.S.-1} is better known in the
following form:

\begin{corollary}
\label{cor.ps2.2.S.-1alt}Let $n\in\mathbb{N}$. Then,%
\[
\sum_{i=0}^{n}\left(  -1\right)  ^{i}\dbinom{n-i}{i}=%
\begin{cases}
1, & \text{if }n\equiv0\operatorname{mod}6\text{ or }n\equiv
1\operatorname{mod}6;\\
0, & \text{if }n\equiv2\operatorname{mod}6\text{ or }n\equiv
5\operatorname{mod}6;\\
-1, & \text{if }n\equiv3\operatorname{mod}6\text{ or }n\equiv
4\operatorname{mod}6
\end{cases}
.
\]

\end{corollary}

\begin{proof}
[Proof of Corollary \ref{cor.ps2.2.S.-1alt}.]Corollary \ref{cor.ps2.2.S.-1}
(applied to $n+1$ instead of $n$) yields%
\[
\sum_{i=0}^{\left(  n+1\right)  -1}\left(  -1\right)  ^{i}\dbinom{\left(
n+1\right)  -1-i}{i}=\left(  -1\right)  ^{n+1}\cdot%
\begin{cases}
0, & \text{if }n+1\equiv0\operatorname{mod}3;\\
-1, & \text{if }n+1\equiv1\operatorname{mod}3;\\
1, & \text{if }n+1\equiv2\operatorname{mod}3
\end{cases}
.
\]
In view of $\left(  n+1\right)  -1=n$, this rewrites as%
\begin{equation}
\sum_{i=0}^{n}\left(  -1\right)  ^{i}\dbinom{n-i}{i}=\left(  -1\right)
^{n+1}\cdot%
\begin{cases}
0, & \text{if }n+1\equiv0\operatorname{mod}3;\\
-1, & \text{if }n+1\equiv1\operatorname{mod}3;\\
1, & \text{if }n+1\equiv2\operatorname{mod}3
\end{cases}
. \label{pf.cor.ps2.2.S.-1alt.1}%
\end{equation}

\begin{vershort}
All that remains to be done now is to prove that the right hand side of
(\ref{pf.cor.ps2.2.S.-1alt.1}) equals
\[%
\begin{cases}
1, & \text{if }n\equiv0\operatorname{mod}6\text{ or }n\equiv
1\operatorname{mod}6;\\
0, & \text{if }n\equiv2\operatorname{mod}6\text{ or }n\equiv
5\operatorname{mod}6;\\
-1, & \text{if }n\equiv3\operatorname{mod}6\text{ or }n\equiv
4\operatorname{mod}6
\end{cases}
.
\]
This can be verified by straightforward computation, distinguishing six
possible cases (one for each remainder $n$ can leave upon division by $6$).
\end{vershort}

\begin{verlong}
The remainder obtained when $n$ is divided by $6$ must be either $0$ or $1$ or
$2$ or $3$ or $4$ or $5$. Hence, we must have either $n\equiv
0\operatorname{mod}6$ or $n\equiv1\operatorname{mod}6$ or $n\equiv
2\operatorname{mod}6$ or $n\equiv3\operatorname{mod}6$ or $n\equiv
4\operatorname{mod}6$ or $n\equiv5\operatorname{mod}6$. Thus, we are in one of
the following six cases:

\textit{Case 1:} We have $n\equiv0\operatorname{mod}6$.

\textit{Case 2:} We have $n\equiv1\operatorname{mod}6$.

\textit{Case 3:} We have $n\equiv2\operatorname{mod}6$.

\textit{Case 4:} We have $n\equiv3\operatorname{mod}6$.

\textit{Case 5:} We have $n\equiv4\operatorname{mod}6$.

\textit{Case 6:} We have $n\equiv5\operatorname{mod}6$.

Let us first consider Case 1. In this case, we have $n\equiv
0\operatorname{mod}6$. Thus, $\underbrace{n}_{\equiv0\operatorname{mod}%
6}+1\equiv0+1=1\operatorname{mod}6$. Therefore, $n+1\equiv1\operatorname{mod}%
2$ and $n+1\equiv1\operatorname{mod}3$. Now, (\ref{pf.cor.ps2.2.S.-1alt.1})
becomes%
\begin{align*}
\sum_{i=0}^{n}\left(  -1\right)  ^{i}\dbinom{n-i}{i}  &  =\underbrace{\left(
-1\right)  ^{n+1}}_{\substack{=-1\\\text{(since }n+1\equiv1\operatorname{mod}%
2\text{)}}}\cdot\underbrace{%
\begin{cases}
0, & \text{if }n+1\equiv0\operatorname{mod}3;\\
-1, & \text{if }n+1\equiv1\operatorname{mod}3;\\
1, & \text{if }n+1\equiv2\operatorname{mod}3
\end{cases}
}_{\substack{=-1\\\text{(since }n+1\equiv1\operatorname{mod}3\text{)}%
}}=\left(  -1\right)  \left(  -1\right) \\
&  =1.
\end{align*}
Comparing this with%
\[%
\begin{cases}
1, & \text{if }n\equiv0\operatorname{mod}6\text{ or }n\equiv
1\operatorname{mod}6;\\
0, & \text{if }n\equiv2\operatorname{mod}6\text{ or }n\equiv
5\operatorname{mod}6;\\
-1, & \text{if }n\equiv3\operatorname{mod}6\text{ or }n\equiv
4\operatorname{mod}6
\end{cases}
=1\ \ \ \ \ \ \ \ \ \ \left(
\begin{array}
[c]{c}%
\text{since }n\equiv0\operatorname{mod}6\text{ or }n\equiv1\operatorname{mod}%
6\\
\text{(because }n\equiv0\operatorname{mod}6\text{)}%
\end{array}
\right)  ,
\]
we obtain%
\[
\sum_{i=0}^{n}\left(  -1\right)  ^{i}\dbinom{n-i}{i}=%
\begin{cases}
1, & \text{if }n\equiv0\operatorname{mod}6\text{ or }n\equiv
1\operatorname{mod}6;\\
0, & \text{if }n\equiv2\operatorname{mod}6\text{ or }n\equiv
5\operatorname{mod}6;\\
-1, & \text{if }n\equiv3\operatorname{mod}6\text{ or }n\equiv
4\operatorname{mod}6
\end{cases}
.
\]
Hence, Corollary \ref{cor.ps2.2.S.-1alt} is proven in Case 1.

Let us next consider Case 2. In this case, we have $n\equiv1\operatorname{mod}%
6$. Thus, $\underbrace{n}_{\equiv1\operatorname{mod}6}+1\equiv
1+1=2\operatorname{mod}6$. Therefore, $n+1\equiv2\equiv0\operatorname{mod}2$
and $n+1\equiv2\operatorname{mod}3$. Now, (\ref{pf.cor.ps2.2.S.-1alt.1})
becomes%
\begin{align*}
\sum_{i=0}^{n}\left(  -1\right)  ^{i}\dbinom{n-i}{i}  &  =\underbrace{\left(
-1\right)  ^{n+1}}_{\substack{=1\\\text{(since }n+1\equiv0\operatorname{mod}%
2\text{)}}}\cdot\underbrace{%
\begin{cases}
0, & \text{if }n+1\equiv0\operatorname{mod}3;\\
-1, & \text{if }n+1\equiv1\operatorname{mod}3;\\
1, & \text{if }n+1\equiv2\operatorname{mod}3
\end{cases}
}_{\substack{=1\\\text{(since }n+1\equiv2\operatorname{mod}3\text{)}}%
}=1\cdot1\\
&  =1.
\end{align*}
Comparing this with%
\[%
\begin{cases}
1, & \text{if }n\equiv0\operatorname{mod}6\text{ or }n\equiv
1\operatorname{mod}6;\\
0, & \text{if }n\equiv2\operatorname{mod}6\text{ or }n\equiv
5\operatorname{mod}6;\\
-1, & \text{if }n\equiv3\operatorname{mod}6\text{ or }n\equiv
4\operatorname{mod}6
\end{cases}
=1\ \ \ \ \ \ \ \ \ \ \left(
\begin{array}
[c]{c}%
\text{since }n\equiv0\operatorname{mod}6\text{ or }n\equiv1\operatorname{mod}%
6\\
\text{(because }n\equiv1\operatorname{mod}6\text{)}%
\end{array}
\right)  ,
\]
we obtain%
\[
\sum_{i=0}^{n}\left(  -1\right)  ^{i}\dbinom{n-i}{i}=%
\begin{cases}
1, & \text{if }n\equiv0\operatorname{mod}6\text{ or }n\equiv
1\operatorname{mod}6;\\
0, & \text{if }n\equiv2\operatorname{mod}6\text{ or }n\equiv
5\operatorname{mod}6;\\
-1, & \text{if }n\equiv3\operatorname{mod}6\text{ or }n\equiv
4\operatorname{mod}6
\end{cases}
.
\]
Hence, Corollary \ref{cor.ps2.2.S.-1alt} is proven in Case 2.

Let us next consider Case 3. In this case, we have $n\equiv2\operatorname{mod}%
6$. Thus, $\underbrace{n}_{\equiv2\operatorname{mod}6}+1\equiv
2+1=3\operatorname{mod}6$. Therefore, $n+1\equiv3\equiv1\operatorname{mod}2$
and $n+1\equiv3\equiv0\operatorname{mod}3$. Now, (\ref{pf.cor.ps2.2.S.-1alt.1}%
) becomes%
\begin{align*}
\sum_{i=0}^{n}\left(  -1\right)  ^{i}\dbinom{n-i}{i}  &  =\left(  -1\right)
^{n+1}\cdot\underbrace{%
\begin{cases}
0, & \text{if }n+1\equiv0\operatorname{mod}3;\\
-1, & \text{if }n+1\equiv1\operatorname{mod}3;\\
1, & \text{if }n+1\equiv2\operatorname{mod}3
\end{cases}
}_{\substack{=0\\\text{(since }n+1\equiv0\operatorname{mod}3\text{)}}}=\left(
-1\right)  ^{n+1}\cdot0\\
&  =0.
\end{align*}
Comparing this with%
\[%
\begin{cases}
1, & \text{if }n\equiv0\operatorname{mod}6\text{ or }n\equiv
1\operatorname{mod}6;\\
0, & \text{if }n\equiv2\operatorname{mod}6\text{ or }n\equiv
5\operatorname{mod}6;\\
-1, & \text{if }n\equiv3\operatorname{mod}6\text{ or }n\equiv
4\operatorname{mod}6
\end{cases}
=1\ \ \ \ \ \ \ \ \ \ \left(
\begin{array}
[c]{c}%
\text{since }n\equiv2\operatorname{mod}6\text{ or }n\equiv5\operatorname{mod}%
6\\
\text{(because }n\equiv2\operatorname{mod}6\text{)}%
\end{array}
\right)  ,
\]
we obtain%
\[
\sum_{i=0}^{n}\left(  -1\right)  ^{i}\dbinom{n-i}{i}=%
\begin{cases}
1, & \text{if }n\equiv0\operatorname{mod}6\text{ or }n\equiv
1\operatorname{mod}6;\\
0, & \text{if }n\equiv2\operatorname{mod}6\text{ or }n\equiv
5\operatorname{mod}6;\\
-1, & \text{if }n\equiv3\operatorname{mod}6\text{ or }n\equiv
4\operatorname{mod}6
\end{cases}
.
\]
Hence, Corollary \ref{cor.ps2.2.S.-1alt} is proven in Case 3.

Let us next consider Case 4. In this case, we have $n\equiv3\operatorname{mod}%
6$. Thus, $\underbrace{n}_{\equiv3\operatorname{mod}6}+1\equiv
3+1=4\operatorname{mod}6$. Therefore, $n+1\equiv4\equiv0\operatorname{mod}2$
and $n+1\equiv4\equiv1\operatorname{mod}3$. Now, (\ref{pf.cor.ps2.2.S.-1alt.1}%
) becomes%
\begin{align*}
\sum_{i=0}^{n}\left(  -1\right)  ^{i}\dbinom{n-i}{i}  &  =\underbrace{\left(
-1\right)  ^{n+1}}_{\substack{=1\\\text{(since }n+1\equiv0\operatorname{mod}%
2\text{)}}}\cdot\underbrace{%
\begin{cases}
0, & \text{if }n+1\equiv0\operatorname{mod}3;\\
-1, & \text{if }n+1\equiv1\operatorname{mod}3;\\
1, & \text{if }n+1\equiv2\operatorname{mod}3
\end{cases}
}_{\substack{=-1\\\text{(since }n+1\equiv1\operatorname{mod}3\text{)}}%
}=1\cdot\left(  -1\right) \\
&  =-1.
\end{align*}
Comparing this with%
\[%
\begin{cases}
1, & \text{if }n\equiv0\operatorname{mod}6\text{ or }n\equiv
1\operatorname{mod}6;\\
0, & \text{if }n\equiv2\operatorname{mod}6\text{ or }n\equiv
5\operatorname{mod}6;\\
-1, & \text{if }n\equiv3\operatorname{mod}6\text{ or }n\equiv
4\operatorname{mod}6
\end{cases}
=-1\ \ \ \ \ \ \ \ \ \ \left(
\begin{array}
[c]{c}%
\text{since }n\equiv3\operatorname{mod}6\text{ or }n\equiv4\operatorname{mod}%
6\\
\text{(because }n\equiv3\operatorname{mod}6\text{)}%
\end{array}
\right)  ,
\]
we obtain%
\[
\sum_{i=0}^{n}\left(  -1\right)  ^{i}\dbinom{n-i}{i}=%
\begin{cases}
1, & \text{if }n\equiv0\operatorname{mod}6\text{ or }n\equiv
1\operatorname{mod}6;\\
0, & \text{if }n\equiv2\operatorname{mod}6\text{ or }n\equiv
5\operatorname{mod}6;\\
-1, & \text{if }n\equiv3\operatorname{mod}6\text{ or }n\equiv
4\operatorname{mod}6
\end{cases}
.
\]
Hence, Corollary \ref{cor.ps2.2.S.-1alt} is proven in Case 4.

Let us next consider Case 5. In this case, we have $n\equiv4\operatorname{mod}%
6$. Thus, $\underbrace{n}_{\equiv4\operatorname{mod}6}+1\equiv
4+1=5\operatorname{mod}6$. Therefore, $n+1\equiv5\equiv1\operatorname{mod}2$
and $n+1\equiv5\equiv2\operatorname{mod}3$. Now, (\ref{pf.cor.ps2.2.S.-1alt.1}%
) becomes%
\begin{align*}
\sum_{i=0}^{n}\left(  -1\right)  ^{i}\dbinom{n-i}{i}  &  =\underbrace{\left(
-1\right)  ^{n+1}}_{\substack{=-1\\\text{(since }n+1\equiv1\operatorname{mod}%
2\text{)}}}\cdot\underbrace{%
\begin{cases}
0, & \text{if }n+1\equiv0\operatorname{mod}3;\\
-1, & \text{if }n+1\equiv1\operatorname{mod}3;\\
1, & \text{if }n+1\equiv2\operatorname{mod}3
\end{cases}
}_{\substack{=1\\\text{(since }n+1\equiv2\operatorname{mod}3\text{)}}}=\left(
-1\right)  \cdot1\\
&  =-1.
\end{align*}
Comparing this with%
\[%
\begin{cases}
1, & \text{if }n\equiv0\operatorname{mod}6\text{ or }n\equiv
1\operatorname{mod}6;\\
0, & \text{if }n\equiv2\operatorname{mod}6\text{ or }n\equiv
5\operatorname{mod}6;\\
-1, & \text{if }n\equiv3\operatorname{mod}6\text{ or }n\equiv
4\operatorname{mod}6
\end{cases}
=1\ \ \ \ \ \ \ \ \ \ \left(
\begin{array}
[c]{c}%
\text{since }n\equiv3\operatorname{mod}6\text{ or }n\equiv4\operatorname{mod}%
6\\
\text{(because }n\equiv4\operatorname{mod}6\text{)}%
\end{array}
\right)  ,
\]
we obtain%
\[
\sum_{i=0}^{n}\left(  -1\right)  ^{i}\dbinom{n-i}{i}=%
\begin{cases}
1, & \text{if }n\equiv0\operatorname{mod}6\text{ or }n\equiv
1\operatorname{mod}6;\\
0, & \text{if }n\equiv2\operatorname{mod}6\text{ or }n\equiv
5\operatorname{mod}6;\\
-1, & \text{if }n\equiv3\operatorname{mod}6\text{ or }n\equiv
4\operatorname{mod}6
\end{cases}
.
\]
Hence, Corollary \ref{cor.ps2.2.S.-1alt} is proven in Case 5.

Let us next consider Case 6. In this case, we have $n\equiv5\operatorname{mod}%
6$. Thus, $\underbrace{n}_{\equiv5\operatorname{mod}6}+1\equiv5+1\equiv
0\operatorname{mod}6$. Therefore, $n+1\equiv0\operatorname{mod}2$ and
$n+1\equiv0\operatorname{mod}3$. Now, (\ref{pf.cor.ps2.2.S.-1alt.1}) becomes%
\begin{align*}
\sum_{i=0}^{n}\left(  -1\right)  ^{i}\dbinom{n-i}{i}  &  =\left(  -1\right)
^{n+1}\cdot\underbrace{%
\begin{cases}
0, & \text{if }n+1\equiv0\operatorname{mod}3;\\
-1, & \text{if }n+1\equiv1\operatorname{mod}3;\\
1, & \text{if }n+1\equiv2\operatorname{mod}3
\end{cases}
}_{\substack{=0\\\text{(since }n+1\equiv0\operatorname{mod}3\text{)}}}=\left(
-1\right)  ^{n+1}\cdot0\\
&  =0.
\end{align*}
Comparing this with%
\[%
\begin{cases}
1, & \text{if }n\equiv0\operatorname{mod}6\text{ or }n\equiv
1\operatorname{mod}6;\\
0, & \text{if }n\equiv2\operatorname{mod}6\text{ or }n\equiv
5\operatorname{mod}6;\\
-1, & \text{if }n\equiv3\operatorname{mod}6\text{ or }n\equiv
4\operatorname{mod}6
\end{cases}
=1\ \ \ \ \ \ \ \ \ \ \left(
\begin{array}
[c]{c}%
\text{since }n\equiv2\operatorname{mod}6\text{ or }n\equiv5\operatorname{mod}%
6\\
\text{(because }n\equiv5\operatorname{mod}6\text{)}%
\end{array}
\right)  ,
\]
we obtain%
\[
\sum_{i=0}^{n}\left(  -1\right)  ^{i}\dbinom{n-i}{i}=%
\begin{cases}
1, & \text{if }n\equiv0\operatorname{mod}6\text{ or }n\equiv
1\operatorname{mod}6;\\
0, & \text{if }n\equiv2\operatorname{mod}6\text{ or }n\equiv
5\operatorname{mod}6;\\
-1, & \text{if }n\equiv3\operatorname{mod}6\text{ or }n\equiv
4\operatorname{mod}6
\end{cases}
.
\]
Hence, Corollary \ref{cor.ps2.2.S.-1alt} is proven in Case 6.

We have now proven Corollary \ref{cor.ps2.2.S.-1alt} in all six Cases 1, 2, 3,
4, 5 and 6. Since these six Cases cover all possibilities, we thus conclude
that Corollary \ref{cor.ps2.2.S.-1alt} always holds.
\end{verlong}

Thus, Corollary \ref{cor.ps2.2.S.-1alt} is proven.
\end{proof}

See also \cite[Identity 172]{BenQui03} or \cite{BenQui08} for a combinatorial
proof of Corollary \ref{cor.ps2.2.S.-1alt}. Also, Corollary
\ref{cor.ps2.2.S.-1alt} appears in \cite[\S 5.2, Problem 3]{GKP}.

\subsection{Solution to Exercise \ref{exe.ps2.2.4}}

\begin{proof}
[Solution to Exercise \ref{exe.ps2.2.4}.]\textbf{(a)} This proof is completely
straightforward, and would be left to the reader in any research paper; we
give a few details only:

Let $i\in\left\{  1,2,\ldots,n-2\right\}  $. We need to prove that $s_{i}\circ
s_{i+1}\circ s_{i}=s_{i+1}\circ s_{i}\circ s_{i+1}$. In order to do so, it is
clearly enough to show that $\left(  s_{i}\circ s_{i+1}\circ s_{i}\right)
\left(  h\right)  =\left(  s_{i+1}\circ s_{i}\circ s_{i+1}\right)  \left(
h\right)  $ for every $h\in\left\{  1,2,\ldots,n\right\}  $. So let us fix
$h\in\left\{  1,2,\ldots,n\right\}  $. We must be in one of the following four cases:

\textit{Case 1:} We have $h=i$.

\textit{Case 2:} We have $h=i+1$.

\textit{Case 3:} We have $h=i+2$.

\textit{Case 4:} We have $h\notin\left\{  i,i+1,i+2\right\}  $.

Let us first consider Case 1. In this case, we have $h=i$ and thus
\begin{align*}
\left(  s_{i}\circ s_{i+1}\circ s_{i}\right)  \left(  \underbrace{h}%
_{=i}\right)   &  =\left(  s_{i}\circ s_{i+1}\circ s_{i}\right)  \left(
i\right)  =s_{i}\left(  s_{i+1}\left(  \underbrace{s_{i}\left(  i\right)
}_{=i+1}\right)  \right) \\
&  =s_{i}\left(  \underbrace{s_{i+1}\left(  i+1\right)  }_{=i+2}\right)
=s_{i}\left(  i+2\right)  =i+2.
\end{align*}
A similar computation shows $\left(  s_{i+1}\circ s_{i}\circ s_{i+1}\right)
\left(  h\right)  =i+2$. Thus, $\left(  s_{i}\circ s_{i+1}\circ s_{i}\right)
\left(  h\right)  =i+2=\left(  s_{i+1}\circ s_{i}\circ s_{i+1}\right)  \left(
h\right)  $. Hence, we have proven the equality $\left(  s_{i}\circ
s_{i+1}\circ s_{i}\right)  \left(  h\right)  =\left(  s_{i+1}\circ s_{i}\circ
s_{i+1}\right)  \left(  h\right)  $ in Case 1.

Similarly, we can prove the same equality in Cases 2 and 3.

Now, let us consider Case 4. In this case, we have $h\notin\left\{
i,i+1,i+2\right\}  $. Thus, $h$ is neither $i$ nor $i+1$, so that we have
$s_{i}\left(  h\right)  =h$. Also, $h$ is neither $i+1$ nor $i+2$ (since
$h\notin\left\{  i,i+1,i+2\right\}  $), and thus we have $s_{i+1}\left(
h\right)  =h$. Hence,
\[
\left(  s_{i}\circ s_{i+1}\circ s_{i}\right)  \left(  h\right)  =s_{i}\left(
s_{i+1}\left(  \underbrace{s_{i}\left(  h\right)  }_{=h}\right)  \right)
=s_{i}\left(  \underbrace{s_{i+1}\left(  h\right)  }_{=h}\right)
=s_{i}\left(  h\right)  =h.
\]
Similarly, $\left(  s_{i+1}\circ s_{i}\circ s_{i+1}\right)  \left(  h\right)
=h$. Thus, $\left(  s_{i}\circ s_{i+1}\circ s_{i}\right)  \left(  h\right)
=h=\left(  s_{i+1}\circ s_{i}\circ s_{i+1}\right)  \left(  h\right)  $. Hence,
we have proven $\left(  s_{i}\circ s_{i+1}\circ s_{i}\right)  \left(
h\right)  =\left(  s_{i+1}\circ s_{i}\circ s_{i+1}\right)  \left(  h\right)  $
in Case 4.

Altogether, we have proven the equality $\left(  s_{i}\circ s_{i+1}\circ
s_{i}\right)  \left(  h\right)  =\left(  s_{i+1}\circ s_{i}\circ
s_{i+1}\right)  \left(  h\right)  $ in each of the four Cases 1, 2, 3, and 4.
Thus, $\left(  s_{i}\circ s_{i+1}\circ s_{i}\right)  \left(  h\right)
=\left(  s_{i+1}\circ s_{i}\circ s_{i+1}\right)  \left(  h\right)  $ always
holds. Exercise \ref{exe.ps2.2.4} \textbf{(a)} is thus solved.

\textbf{(b)} We follow the hint and delay the solution of this part until
later (see Exercise \ref{exe.ps2.2.5} \textbf{(e)}).

\textbf{(c)} For every $i\in\left\{  1,2,\ldots,n\right\}  $, we let $a_{i}$
be the permutation $s_{i-1}\circ s_{i-2}\circ\cdots\circ s_{1}\in S_{n}%
$\ \ \ \ \footnote{In particular, $a_{1}=s_{1-1}\circ s_{1-2}\circ\cdots\circ
s_{1}$ is the composition of $0$ permutations. What does this mean? Just as a
sum of $0$ terms is defined to be $0$ (because $0$ is the neutral element of
addition), and a product of $0$ terms is $1$ (since $1$ is the neutral element
of multiplication), the composition of $0$ permutations is defined to be the
identity permutation (since the identity permutation is the neutral element of
composition). Hence, $a_{1}$ is the identity permutation, i.e., we have
$a_{1}=\operatorname*{id}$.}. Now, we claim that%
\begin{equation}
w_{0}=a_{1}\circ a_{2}\circ\cdots\circ a_{n}. \label{sol.ps2.2.4.c.claim}%
\end{equation}
This is essentially an explicit way to write $w_{0}$ as a composition of
several permutations of the form $s_{i}$ (because each $a_{i}$ on the right
hand side is the composition $s_{i-1}\circ s_{i-2}\circ\cdots\circ s_{1}$).
Thus, once (\ref{sol.ps2.2.4.c.claim}) is proven, the exercise will be solved.

Before we prove (\ref{sol.ps2.2.4.c.claim}), let us first understand how
$a_{i}$ operates. We claim that%
\begin{equation}
a_{i}\left(  k\right)  =
\begin{cases}
i, & \text{if }k=1;\\
k-1, & \text{if }1<k\leq i;\\
k, & \text{if }k>i
\end{cases}
\label{sol.ps2.2.4.c.aik}%
\end{equation}
for each $i\in\left\{  1,2,\ldots,n\right\}  $ and $k\in\left\{
1,2,\ldots,n\right\}  $. (In other words, $a_{i}$ is the permutation which
sends $1,2,3,\ldots,i$ to $i,1,2,\ldots,i-1$ and leaves all numbers $>i$
untouched. Using the terminology of Definition \ref{def.perm.cycles}, this
means that $a_{i}=\operatorname*{cyc}\nolimits_{i,i-1,\ldots,1}$.)

[\textit{Proof of (\ref{sol.ps2.2.4.c.aik}):} Let me first give an informal
proof of (\ref{sol.ps2.2.4.c.aik}) which is not hard to turn into a formal proof.

Let $i\in\left\{  1,2,\ldots,n\right\}  $. We have $a_{i}=s_{i-1}\circ
s_{i-2}\circ\cdots\circ s_{1}$. Therefore, $a_{i}$ is the permutation that
first swaps $1$ with $2$, then swaps $2$ with $3$, etc., until it finally
swaps $i-1$ with $i$. Thus:

\begin{itemize}
\item When we apply $a_{i}$ to $1$, we arrive at $i$ at the end (since the $1$
is carried to $2$ by the first swap, which is then carried to $3$ by the next
swap, and so on, until it finally becomes $i$).

\item When we apply $a_{i}$ to some $k\in\left\{  2,3,\ldots,i\right\}  $, we
arrive at $k-1$ at the end (since the first swap to move $k$ is the $\left(
k-1\right)  $-st swap, which changes it into $k-1$, and from then on all the
following swaps leave $k-1$ untouched).

\item When we apply $a_{i}$ to some $k\in\left\{  i+1,i+2,\ldots,n\right\}  $,
we arrive at $k$ at the end (since none of our swaps changes $k$).
\end{itemize}

Expressing this in a formula instead of in words, we obtain precisely
(\ref{sol.ps2.2.4.c.aik}).

\textit{Formal proof of (\ref{sol.ps2.2.4.c.aik}):} For the sake of
completeness, let me show how to prove (\ref{sol.ps2.2.4.c.aik}) formally.

We shall prove (\ref{sol.ps2.2.4.c.aik}) by induction on $i$:

\textit{Induction base:} We have $a_{1}=s_{0}\circ s_{-1}\circ\cdots\circ
s_{1}=\left(  \text{a composition of }0\text{ permutations}\right)
=\operatorname*{id}$. Thus, every $k\in\left\{  1,2,\ldots,n\right\}  $
satisfies%
\begin{align*}
\underbrace{a_{1}}_{=\operatorname*{id}}\left(  k\right)   &
=\operatorname*{id}\left(  k\right)  =k =
\begin{cases}
k, & \text{if }k=1;\\
k, & \text{if }k>1
\end{cases}
=
\begin{cases}
1, & \text{if }k=1;\\
k, & \text{if }k>1
\end{cases}
\\
&  \ \ \ \ \ \ \ \ \ \ \left(  \text{since }k=1\text{ when }k=1\right) \\
&  =
\begin{cases}
1, & \text{if }k=1;\\
k-1, & \text{if }1<k\leq1;\\
k, & \text{if }k>1
\end{cases}
\\
&  \ \ \ \ \ \ \ \ \ \ \left(
\begin{array}
[c]{c}%
\text{here, we added a \textquotedblleft}1<k\leq1\text{\textquotedblright%
\ case, which does not}\\
\text{change the result because this case never happens}%
\end{array}
\right)  .
\end{align*}
In other words, (\ref{sol.ps2.2.4.c.aik}) holds for $i=1$. This completes the
induction base.

\textit{Induction step:} Let $I\in\left\{  1,2,\ldots,n\right\}  $ be such
that $I>1$. Assume that (\ref{sol.ps2.2.4.c.aik}) holds for $i=I-1$. We need
to show that (\ref{sol.ps2.2.4.c.aik}) holds for $i=I$.

We have assumed that (\ref{sol.ps2.2.4.c.aik}) holds for $i=I-1$. In other
words,%
\begin{equation}
a_{I-1}\left(  k\right)  =
\begin{cases}
I-1, & \text{if }k=1;\\
k-1, & \text{if }1<k\leq I-1;\\
k, & \text{if }k>I-1
\end{cases}
\label{sol.ps2.2.4.c.aik.pf.hyp}%
\end{equation}
for each $k\in\left\{  1,2,\ldots,n\right\}  $.

The definition of $a_{I-1}$ yields $a_{I-1}=s_{I-2}\circ s_{I-3}\circ
\cdots\circ s_{1}$. The definition of $a_{I}$ yields $a_{I}=s_{I-1}\circ
s_{I-2}\circ\cdots\circ s_{1}=s_{I-1}\circ\underbrace{\left(  s_{I-2}\circ
s_{I-3}\circ\cdots\circ s_{1}\right)  }_{=a_{I-1}}=s_{I-1}\circ a_{I-1}$.
Hence, every $k\in\left\{  1,2,\ldots,n\right\}  $ satisfies
\begin{equation}
a_{I}\left(  k\right)  =
\begin{cases}
I, & \text{if }k=1;\\
k-1, & \text{if }1<k\leq I;\\
k, & \text{if }k>I
\end{cases}
\label{sol.ps2.2.4.c.aik.pf.step}%
\end{equation}
\footnote{\textit{Proof of (\ref{sol.ps2.2.4.c.aik.pf.step}):} Let
$k\in\left\{  1,2,\ldots,n\right\}  $. We need to prove
(\ref{sol.ps2.2.4.c.aik.pf.step}). We are in one of the following four cases:
\par
\textit{Case 1:} We have $k=1$.
\par
\textit{Case 2:} We have $1<k\leq I$ and $k<I$.
\par
\textit{Case 3:} We have $1<k\leq I$ and $k\geq I$.
\par
\textit{Case 4:} We have $k>I$.
\par
Let us first consider Case 1. In this case, we have $k=1$. Thus,
(\ref{sol.ps2.2.4.c.aik.pf.hyp}) yields $a_{I-1}\left(  k\right)  =
\begin{cases}
I-1, & \text{if }k=1;\\
k-1, & \text{if }1<k\leq I-1;\\
k, & \text{if }k>I-1
\end{cases}
=I-1$ (since $k=1$) and thus $\underbrace{a_{I}}_{=s_{I-1}\circ a_{I-1}%
}\left(  k\right)  =\left(  s_{I-1}\circ a_{I-1}\right)  \left(  k\right)
=s_{I-1}\left(  \underbrace{a_{I-1}\left(  k\right)  }_{=I-1}\right)
=s_{I-1}\left(  I-1\right)  =I$. Compared with $%
\begin{cases}
I, & \text{if }k=1;\\
k-1, & \text{if }1<k\leq I;\\
k, & \text{if }k>I
\end{cases}
=I$ (since $k=1$), this yields $a_{I}\left(  k\right)  =
\begin{cases}
I, & \text{if }k=1;\\
k-1, & \text{if }1<k\leq I;\\
k, & \text{if }k>I
\end{cases}
$. Thus, (\ref{sol.ps2.2.4.c.aik.pf.step}) is proven in Case 1.
\par
Let us now consider Case 2. In this case, we have $1<k\leq I$ and $k<I$. From
$k<I$, we obtain $k\leq I-1$ (since $k$ and $I$ are integers). Thus,
(\ref{sol.ps2.2.4.c.aik.pf.hyp}) yields $a_{I-1}\left(  k\right)  =
\begin{cases}
I-1, & \text{if }k=1;\\
k-1, & \text{if }1<k\leq I-1;\\
k, & \text{if }k>I-1
\end{cases}
=k-1$ (since $1<k\leq I-1$) and thus $\underbrace{a_{I}}_{=s_{I-1}\circ
a_{I-1}}\left(  k\right)  =\left(  s_{I-1}\circ a_{I-1}\right)  \left(
k\right)  =s_{I-1}\left(  \underbrace{a_{I-1}\left(  k\right)  }%
_{=k-1}\right)  =s_{I-1}\left(  k-1\right)  =k-1$ (since $k-1<k\leq I-1$).
Compared with $%
\begin{cases}
I, & \text{if }k=1;\\
k-1, & \text{if }1<k\leq I;\\
k, & \text{if }k>I
\end{cases}
=k-1$ (since $1<k\leq I$), this yields $a_{I}\left(  k\right)  =
\begin{cases}
I, & \text{if }k=1;\\
k-1, & \text{if }1<k\leq I;\\
k, & \text{if }k>I
\end{cases}
$. Thus, (\ref{sol.ps2.2.4.c.aik.pf.step}) is proven in Case 2.
\par
Let us now consider Case 3. In this case, we have $1<k\leq I$ and $k\geq I$.
Combining $k\leq I$ with $k\geq I$, we obtain $k=I>I-1$. Thus,
(\ref{sol.ps2.2.4.c.aik.pf.hyp}) yields $a_{I-1}\left(  k\right)  =
\begin{cases}
I-1, & \text{if }k=1;\\
k-1, & \text{if }1<k\leq I-1;\\
k, & \text{if }k>I-1
\end{cases}
=k$ (since $k>I-1$) and thus $\underbrace{a_{I}}_{=s_{I-1}\circ a_{I-1}%
}\left(  k\right)  =\left(  s_{I-1}\circ a_{I-1}\right)  \left(  k\right)
=s_{I-1}\left(  \underbrace{a_{I-1}\left(  k\right)  }_{=k=I}\right)
=s_{I-1}\left(  I\right)  =\underbrace{I}_{=k}-1=k-1$. Compared with $%
\begin{cases}
I, & \text{if }k=1;\\
k-1, & \text{if }1<k\leq I;\\
k, & \text{if }k>I
\end{cases}
=k-1$ (since $1<k\leq I$), this yields $a_{I}\left(  k\right)  =%
\begin{cases}
I, & \text{if }k=1;\\
k-1, & \text{if }1<k\leq I;\\
k, & \text{if }k>I
\end{cases}
$. Thus, (\ref{sol.ps2.2.4.c.aik.pf.step}) is proven in Case 3.
\par
Let us finally consider Case 4. In this case, we have $k>I$. Hence, $k>I>I-1$.
Thus, (\ref{sol.ps2.2.4.c.aik.pf.hyp}) yields $a_{I-1}\left(  k\right)  =%
\begin{cases}
I-1, & \text{if }k=1;\\
k-1, & \text{if }1<k\leq I-1;\\
k, & \text{if }k>I-1
\end{cases}
=k$ (since $k>I-1$) and thus $\underbrace{a_{I}}_{=s_{I-1}\circ a_{I-1}%
}\left(  k\right)  =\left(  s_{I-1}\circ a_{I-1}\right)  \left(  k\right)
=s_{I-1}\left(  \underbrace{a_{I-1}\left(  k\right)  }_{=k}\right)
=s_{I-1}\left(  k\right)  =k$ (since $k>I$). Compared with $%
\begin{cases}
I, & \text{if }k=1;\\
k-1, & \text{if }1<k\leq I;\\
k, & \text{if }k>I
\end{cases}
=k$ (since $k>I$), this yields $a_{I}\left(  k\right)  =
\begin{cases}
I, & \text{if }k=1;\\
k-1, & \text{if }1<k\leq I;\\
k, & \text{if }k>I
\end{cases}
$. Thus, (\ref{sol.ps2.2.4.c.aik.pf.step}) is proven in Case 4.
\par
Hence, (\ref{sol.ps2.2.4.c.aik.pf.step}) is proven in each of the four Cases
1, 2, 3 and 4. Since these four Cases are the only possible cases, this shows
that (\ref{sol.ps2.2.4.c.aik.pf.step}) always holds, qed.}. In other words,
(\ref{sol.ps2.2.4.c.aik}) holds for $i=I$. This completes the induction step.
The induction proof of (\ref{sol.ps2.2.4.c.aik}) is thus finished.]

Next, for every $m\in\left\{  0,1,\ldots,n\right\}  $, set%
\[
b_{m}=a_{1}\circ a_{2}\circ\cdots\circ a_{m}\in S_{n}.
\]
As a consequence, $b_{0}=a_{1}\circ a_{2}\circ\cdots\circ a_{0}=\left(
\text{a composition of }0\text{ maps}\right)  =\operatorname*{id}$ and
$b_{n}=a_{1}\circ a_{2}\circ\cdots\circ a_{n}$. We claim that%
\begin{equation}
b_{m}\left(  k\right)  =
\begin{cases}
m+1-k, & \text{if }k\leq m;\\
k, & \text{if }k>m
\end{cases}
\label{sol.ps2.2.4.c.bik}%
\end{equation}
for every $m\in\left\{  0,1,\ldots,n\right\}  $ and $k\in\left\{
1,2,\ldots,n\right\}  $.

[Again, one can prove (\ref{sol.ps2.2.4.c.bik}) either formally by induction
on $m$, or more intuitively by tracking what happens to $k$ under the maps
$a_{m}$, $a_{m-1}$, $\ldots$, $a_{1}$ when these maps are applied one after
the other. This time the informal way is a bit tricky, so let us show the
formal one in all its glory. (You do not ever need to write a proof in this
level of detail unless you are talking to a computer.)

\textit{Proof of (\ref{sol.ps2.2.4.c.bik}):} We shall prove
(\ref{sol.ps2.2.4.c.bik}) by induction on $m$:

\textit{Induction base:} For every $k\in\left\{  1,2,\ldots,n\right\}  $, we
have $\underbrace{b_{0}}_{=\operatorname*{id}}\left(  k\right)
=\operatorname*{id}\left(  k\right)  =k=
\begin{cases}
0+1-k, & \text{if }k\leq0;\\
k, & \text{if }k>0
\end{cases}
$ (since $%
\begin{cases}
0+1-k, & \text{if }k\leq0;\\
k, & \text{if }k>0
\end{cases}
=k$ (since $k>0$)). In other words, (\ref{sol.ps2.2.4.c.bik}) holds for $m=0$.
The induction base is thus complete.

\textit{Induction step:} Let $M\in\left\{  0,1,\ldots,n\right\}  $ be such
that $M>0$. Assume that (\ref{sol.ps2.2.4.c.bik}) holds for $m=M-1$. We need
to show that (\ref{sol.ps2.2.4.c.bik}) holds for $m=M$.

We have assumed that (\ref{sol.ps2.2.4.c.bik}) holds for $m=M-1$. In other
words,
\[
b_{M-1}\left(  k\right)  =
\begin{cases}
\left(  M-1\right)  +1-k, & \text{if }k\leq M-1;\\
k, & \text{if }k>M-1
\end{cases}
\]
for every $k\in\left\{  1,2,\ldots,n\right\}  $. Thus, for every $k\in\left\{
1,2,\ldots,n\right\}  $, we have%
\begin{align}
b_{M-1}\left(  k\right)   &  =
\begin{cases}
\left(  M-1\right)  +1-k, & \text{if }k\leq M-1;\\
k, & \text{if }k>M-1
\end{cases}
\nonumber\\
&  =%
\begin{cases}
M-k, & \text{if }k\leq M-1;\\
k, & \text{if }k>M-1
\end{cases}
\label{sol.ps2.2.4.c.bik.pf.hyp}%
\end{align}
(since $\left(  M-1\right)  +1-k=M-k$).

The definition of $b_{M-1}$ yields $b_{M-1}=a_{1}\circ a_{2}\circ\cdots\circ
a_{M-1}$. The definition of $b_{M}$ yields $b_{M}=a_{1}\circ a_{2}\circ
\cdots\circ a_{M}=\underbrace{\left(  a_{1}\circ a_{2}\circ\cdots\circ
a_{M-1}\right)  }_{=b_{M-1}}\circ a_{M}=b_{M-1}\circ a_{M}$. Thus, we obtain%
\begin{equation}
b_{M}\left(  k\right)  =
\begin{cases}
M+1-k, & \text{if }k\leq M;\\
k, & \text{if }k>M
\end{cases}
\label{sol.ps2.2.4.c.bik.pf.step}%
\end{equation}
for every $k\in\left\{  1,2,\ldots,n\right\}  $\ \ \ \ \footnote{\textit{Proof
of (\ref{sol.ps2.2.4.c.bik.pf.step}):} Let $k\in\left\{  1,2,\ldots,n\right\}
$. We need to prove (\ref{sol.ps2.2.4.c.bik.pf.step}). We are in one of the
following three cases:
\par
\textit{Case 1:} We have $k=1$.
\par
\textit{Case 2:} We have $1<k\leq M$.
\par
\textit{Case 3:} We have $k>M$.
\par
Let us first consider Case 1. In this case, we have $k=1$. Thus, $k=1\leq M$
(since $M\geq1$ (since $M>0$)). But (\ref{sol.ps2.2.4.c.aik}) (applied to
$i=M$) yields $a_{M}\left(  k\right)  =
\begin{cases}
M, & \text{if }k=1;\\
k-1, & \text{if }1<k\leq M;\\
k, & \text{if }k>M
\end{cases}
=M$ (since $k=1$), so that
\begin{align*}
\underbrace{b_{M}}_{=b_{M-1}\circ a_{M}}\left(  k\right)   &  =\left(
b_{M-1}\circ a_{M}\right)  \left(  k\right)  =b_{M-1}\left(  \underbrace{a_{M}%
\left(  k\right)  }_{=M}\right)  =b_{M-1}\left(  M\right) \\
&  =
\begin{cases}
M-M, & \text{if }M\leq M-1;\\
M, & \text{if }M>M-1
\end{cases}
\ \ \ \ \ \ \ \ \ \ \left(  \text{by (\ref{sol.ps2.2.4.c.bik.pf.hyp}), applied
to }M\text{ instead of }k\right) \\
&  =M\ \ \ \ \ \ \ \ \ \ \left(  \text{since }M>M-1\right) \\
&  =M+1-\underbrace{1}_{=k}=M+1-k=
\begin{cases}
M+1-k, & \text{if }k\leq M;\\
k, & \text{if }k>M
\end{cases}
\\
&  \ \ \ \ \ \ \ \ \ \ \left(  \text{since }
\begin{cases}
M+1-k, & \text{if }k\leq M;\\
k, & \text{if }k>M
\end{cases}
=M+1-k\text{ (because }k\leq M\text{)}\right)  .
\end{align*}
Thus, (\ref{sol.ps2.2.4.c.bik.pf.step}) is proven in Case 1.
\par
Let us now consider Case 2. In this case, we have $1<k\leq M$. Now,
(\ref{sol.ps2.2.4.c.aik}) (applied to $i=M$) yields $a_{M}\left(  k\right)  =
\begin{cases}
M, & \text{if }k=1;\\
k-1, & \text{if }1<k\leq M;\\
k, & \text{if }k>M
\end{cases}
=k-1$ (since $1<k\leq M$), so that%
\begin{align*}
\underbrace{b_{M}}_{=b_{M-1}\circ a_{M}}\left(  k\right)   &  =\left(
b_{M-1}\circ a_{M}\right)  \left(  k\right)  =b_{M-1}\left(  \underbrace{a_{M}%
\left(  k\right)  }_{=k-1}\right)  =b_{M-1}\left(  k-1\right) \\
&  =
\begin{cases}
M-\left(  k-1\right)  , & \text{if }k-1\leq M-1;\\
k-1, & \text{if }k-1>M-1
\end{cases}
\ \ \ \ \ \ \ \ \ \ \left(
\begin{array}
[c]{c}%
\text{by (\ref{sol.ps2.2.4.c.bik.pf.hyp}), applied to }k-1\\
\text{instead of }k
\end{array}
\right) \\
&  =M-\left(  k-1\right)  \ \ \ \ \ \ \ \ \ \ \left(  \text{since }k-1\leq
M-1\text{ (since }k\leq M\text{)}\right) \\
&  =M+1-k=
\begin{cases}
M+1-k, & \text{if }k\leq M;\\
k, & \text{if }k>M
\end{cases}
\\
&  \ \ \ \ \ \ \ \ \ \ \left(  \text{since }
\begin{cases}
M+1-k, & \text{if }k\leq M;\\
k, & \text{if }k>M
\end{cases}
=M+1-k\text{ (because }k\leq M\text{)}\right)  .
\end{align*}
Thus, (\ref{sol.ps2.2.4.c.bik.pf.step}) is proven in Case 2.
\par
Let us finally consider Case 3. In this case, we have $k>M$. Now,
(\ref{sol.ps2.2.4.c.aik}) (applied to $i=M$) yields $a_{M}\left(  k\right)  =
\begin{cases}
M, & \text{if }k=1;\\
k-1, & \text{if }1<k\leq M;\\
k, & \text{if }k>M
\end{cases}
=k$ (since $k>M$), so that%
\begin{align*}
\underbrace{b_{M}}_{=b_{M-1}\circ a_{M}}\left(  k\right)   &  =\left(
b_{M-1}\circ a_{M}\right)  \left(  k\right)  =b_{M-1}\left(  \underbrace{a_{M}%
\left(  k\right)  }_{=k}\right)  =b_{M-1}\left(  k\right) \\
&  =
\begin{cases}
M-k, & \text{if }k\leq M-1;\\
k, & \text{if }k>M-1
\end{cases}
\ \ \ \ \ \ \ \ \ \ \left(  \text{by (\ref{sol.ps2.2.4.c.bik.pf.hyp})}\right)
\\
&  =k\ \ \ \ \ \ \ \ \ \ \left(  \text{since }k>M>M-1\right) \\
&  =
\begin{cases}
M+1-k, & \text{if }k\leq M;\\
k, & \text{if }k>M
\end{cases}
\\
&  \ \ \ \ \ \ \ \ \ \ \left(  \text{since }
\begin{cases}
M+1-k, & \text{if }k\leq M;\\
k, & \text{if }k>M
\end{cases}
=k\text{ (because }k>M\text{)}\right)  .
\end{align*}
Thus, (\ref{sol.ps2.2.4.c.bik.pf.step}) is proven in Case 3.
\par
Hence, (\ref{sol.ps2.2.4.c.bik.pf.step}) is proven in each of the three Cases
1, 2 and 3. Since these three Cases are the only possible cases, this shows
that (\ref{sol.ps2.2.4.c.bik.pf.step}) always holds, qed.}. In other words,
(\ref{sol.ps2.2.4.c.bik}) holds for $m=M$. This completes the induction step.
Thus, the induction proof of (\ref{sol.ps2.2.4.c.bik}) is complete.]

Now, every $k\in\left\{  1,2,\ldots,n\right\}  $ satisfies
\begin{align*}
b_{n}\left(  k\right)   &  =
\begin{cases}
n+1-k, & \text{if }k\leq n;\\
k, & \text{if }k>n
\end{cases}
\ \ \ \ \ \ \ \ \ \ \left(  \text{by (\ref{sol.ps2.2.4.c.bik}), applied to
}m=n\right) \\
&  =n+1-k\ \ \ \ \ \ \ \ \ \ \left(  \text{since }k\leq n\right) \\
&  =w_{0}\left(  k\right)  \ \ \ \ \ \ \ \ \ \ \left(  \text{since }%
w_{0}\left(  k\right)  =n+1-k\text{ (by the definition of }w_{0}%
\text{)}\right)  .
\end{align*}
In other words, $b_{n}=w_{0}$, so that $w_{0}=b_{n}=a_{1}\circ a_{2}%
\circ\cdots\circ a_{n}$. This proves (\ref{sol.ps2.2.4.c.claim}). As we know,
this solves Exercise \ref{exe.ps2.2.4} \textbf{(c)}.
\end{proof}

\subsection{Solution to Exercise \ref{exe.ps2.2.5}}

Let us first state Exercise \ref{exe.ps2.2.5} \textbf{(d)} as a separate result:

\begin{proposition}
\label{prop.sol.ps2.2.5.d}Let $n\in\mathbb{N}$. Let $\sigma\in S_{n}$ be a
permutation satisfying $\sigma\left(  1\right)  \leq\sigma\left(  2\right)
\leq\cdots\leq\sigma\left(  n\right)  $. Then, $\sigma=\operatorname*{id}$.
\end{proposition}

Proposition \ref{prop.sol.ps2.2.5.d} essentially says that the only way to
list the numbers $1,2,\ldots,n$ in increasing order is $\left(  1,2,\ldots
,n\right)  $. If you think this is intuitively obvious, you are right. Let me
nevertheless give two proofs (the second of which is completely formal):

\begin{proof}
[First proof of Proposition \ref{prop.sol.ps2.2.5.d}.]First of all, every
$k\in\left\{  1,2,\ldots,n-1\right\}  $ satisfies $\sigma\left(  k\right)
<\sigma\left(  k+1\right)  $\ \ \ \ \footnote{\textit{Proof.} Let
$k\in\left\{  1,2,\ldots,n-1\right\}  $. Then, $\sigma$ is a permutation, thus
injective. Hence, $\sigma\left(  k\right)  \neq\sigma\left(  k+1\right)  $
(since $k\neq k+1$). Combined with $\sigma\left(  k\right)  \leq\sigma\left(
k+1\right)  $ (since $\sigma\left(  1\right)  \leq\sigma\left(  2\right)
\leq\cdots\leq\sigma\left(  n\right)  $), this yields $\sigma\left(  k\right)
<\sigma\left(  k+1\right)  $, qed.}. In other words, $\sigma\left(  1\right)
<\sigma\left(  2\right)  <\cdots<\sigma\left(  n\right)  $.

Let $i\in\left\{  1,2,\ldots,n\right\}  $. Then, $\sigma\left(  i\right)  $ is
the $i$-th smallest among the numbers \newline$\sigma\left(  1\right)
,\sigma\left(  2\right)  ,\ldots,\sigma\left(  n\right)  $ (because
$\sigma\left(  1\right)  <\sigma\left(  2\right)  <\cdots<\sigma\left(
n\right)  $). But since the numbers $\sigma\left(  1\right)  ,\sigma\left(
2\right)  ,\ldots,\sigma\left(  n\right)  $ are just the numbers
$1,2,\ldots,n$ (possibly in a different order)\footnote{since $\sigma$ is a
permutation of $\left\{  1,2,\ldots,n\right\}  $}, it is clear that the $i$-th
smallest among these numbers is $i$. Thus, $\sigma\left(  i\right)  =i$ (since
$\sigma\left(  i\right)  $ is the $i$-th smallest among the numbers
$\sigma\left(  1\right)  ,\sigma\left(  2\right)  ,\ldots,\sigma\left(
n\right)  $). Hence, $\sigma\left(  i\right)  =i=\operatorname*{id}\left(
i\right)  $.

Let us now forget that we fixed $i$. Thus, we have shown that $\sigma\left(
i\right)  =\operatorname*{id}\left(  i\right)  $ for every $i\in\left\{
1,2,\ldots,n\right\}  $. In other words, $\sigma=\operatorname*{id}$.
Proposition \ref{prop.sol.ps2.2.5.d} is thus proven.
\end{proof}

\begin{proof}
[Second proof of Proposition \ref{prop.sol.ps2.2.5.d}.]We shall show that
\begin{equation}
\sigma\left(  i\right)  =i\ \ \ \ \ \ \ \ \ \ \text{for every }i\in\left\{
1,2,\ldots,n\right\}  . \label{sol.ps2.2.5.d.1}%
\end{equation}

[\textit{Proof of (\ref{sol.ps2.2.5.d.1}):} We shall prove
(\ref{sol.ps2.2.5.d.1}) by strong induction over $i$. Thus, we fix some
$I\in\left\{  1,2,\ldots,n\right\}  $, and we assume that
(\ref{sol.ps2.2.5.d.1}) is proven for every $i<I$. We then have to prove that
(\ref{sol.ps2.2.5.d.1}) holds for $i=I$.

We have assumed that (\ref{sol.ps2.2.5.d.1}) is proven for every $i<I$. In
other words,%
\begin{equation}
\sigma\left(  i\right)  =i\ \ \ \ \ \ \ \ \ \ \text{for every }i\in\left\{
1,2,\ldots,n\right\}  \text{ satisfying }i<I. \label{sol.ps2.2.5.d.1.pf.1}%
\end{equation}

We assume (for the sake of contradiction) that $\sigma\left(  I\right)  \neq
I$. But $\sigma$ is a permutation, and thus injective. Hence, from
$\sigma\left(  I\right)  \neq I$, we obtain $\sigma\left(  \sigma\left(
I\right)  \right)  \neq\sigma\left(  I\right)  $. But if $\sigma\left(
I\right)  <I$, then $\sigma\left(  \sigma\left(  I\right)  \right)
=\sigma\left(  I\right)  $ (by (\ref{sol.ps2.2.5.d.1.pf.1}), applied to
$i=\sigma\left(  I\right)  $), which contradicts $\sigma\left(  \sigma\left(
I\right)  \right)  \neq\sigma\left(  I\right)  $. Hence, we cannot have
$\sigma\left(  I\right)  <I$. Thus, we have $\sigma\left(  I\right)  \geq I$.
Combined with $\sigma\left(  I\right)  \neq I$, this yields $\sigma\left(
I\right)  >I$.

Now, let $K=\sigma^{-1}\left(  I\right)  $. Then, $I=\sigma\left(  K\right)
$, so that $\sigma\left(  K\right)  =I\neq\sigma\left(  I\right)  $ and
therefore $K\neq I$. If $K<I$, then $\sigma\left(  K\right)  =K$ (by
(\ref{sol.ps2.2.5.d.1.pf.1}), applied to $i=K$), which contradicts
$\sigma\left(  K\right)  =I\neq K$. Hence, we cannot have $K<I$. We thus have
$I\leq K$.

Now, recall that $\sigma\left(  1\right)  \leq\sigma\left(  2\right)
\leq\cdots\leq\sigma\left(  n\right)  $. In other words, $\sigma\left(
a\right)  \leq\sigma\left(  b\right)  $ for every two elements $a$ and $b$ of
$\left\{  1,2,\ldots,n\right\}  $ satisfying $a\leq b$. Applying this to $a=I$
and $b=K$, we obtain $\sigma\left(  I\right)  \leq\sigma\left(  K\right)  $.
This contradicts $\sigma\left(  I\right)  >I=\sigma\left(  K\right)  $. This
contradiction proves that our assumption (that $\sigma\left(  I\right)  \neq
I$) was wrong. Hence, we must have $\sigma\left(  I\right)  =I$. In other
words, (\ref{sol.ps2.2.5.d.1}) holds for $i=I$. This completes our inductive
proof of (\ref{sol.ps2.2.5.d.1}).]

Now, (\ref{sol.ps2.2.5.d.1}) shows that every $i\in\left\{  1,2,\ldots
,n\right\}  $ satisfies $\sigma\left(  i\right)  =i=\operatorname*{id}\left(
i\right)  $. In other words, $\sigma=\operatorname*{id}$. Proposition
\ref{prop.sol.ps2.2.5.d} is proved again.
\end{proof}

For future use, let us record an easy consequence of Proposition
\ref{prop.sol.ps2.2.5.d}:

\begin{corollary}
\label{cor.sol.ps2.2.5.d2}Let $n\in\mathbb{N}$. Let $\sigma\in S_{n}$ be a
permutation satisfying $\ell\left(  \sigma\right)  =0$. Then, $\sigma
=\operatorname*{id}$.
\end{corollary}

\begin{proof}
[Proof of Corollary \ref{cor.sol.ps2.2.5.d2}.]Let $k\in\left\{  1,2,\ldots
,n-1\right\}  $.

\begin{verlong}
Thus, $1\leq k\leq n-1$, so that $k+1\leq n$ (since $k\leq n-1$).
\end{verlong}

Assume (for the sake of contradiction) that $\sigma\left(  k\right)
>\sigma\left(  k+1\right)  $. Then, $\left(  k,k+1\right)  $ is a pair of
integers satisfying $1\leq k<k+1\leq n$ and $\sigma\left(  k\right)
>\sigma\left(  k+1\right)  $. In other words, $\left(  k,k+1\right)  $ is a
pair $\left(  i,j\right)  $ of integers satisfying $1\leq i<j\leq n$ and
$\sigma\left(  i\right)  >\sigma\left(  j\right)  $. In other words, $\left(
k,k+1\right)  $ is an inversion of $\sigma$ (by the definition of an
\textquotedblleft inversion\textquotedblright). Thus, the permutation $\sigma$
has at least one inversion (namely, $\left(  k,k+1\right)  $).

But the number of inversions of $\sigma$ is $\ell\left(  \sigma\right)  =0$;
in other words, $\sigma$ has no inversions. This contradicts the fact that
$\sigma$ has at least one inversion. This contradiction proves that our
assumption (that $\sigma\left(  k\right)  >\sigma\left(  k+1\right)  $) was
wrong. Hence, we have $\sigma\left(  k\right)  \leq\sigma\left(  k+1\right)  $.

Now, let us forget that we fixed $k$. We thus have shown that $\sigma\left(
k\right)  \leq\sigma\left(  k+1\right)  $ for every $k\in\left\{
1,2,\ldots,n-1\right\}  $. Thus, $\sigma\left(  1\right)  \leq\sigma\left(
2\right)  \leq\cdots\leq\sigma\left(  n\right)  $. Therefore, Proposition
\ref{prop.sol.ps2.2.5.d} shows that $\sigma=\operatorname*{id}$. This proves
Corollary \ref{cor.sol.ps2.2.5.d2}.
\end{proof}

Now, we come to the actual solution of Exercise \ref{exe.ps2.2.5}.

\begin{proof}
[Solution to Exercise \ref{exe.ps2.2.5}.]Exercise \ref{exe.ps2.2.5}
\textbf{(d)} follows immediately from Proposition \ref{prop.sol.ps2.2.5.d}. We
shall next prove part \textbf{(f)} of the exercise, then part \textbf{(a)},
then part \textbf{(e)}, and then the remaining three parts.

Before we come to the actual solution, let us introduce one more notation.

For every $\sigma\in S_{n}$, let $\operatorname*{Inv}\left(  \sigma\right)  $
be the set of all inversions of the permutation $\sigma$. Thus, for every
$\sigma\in S_{n}$, we have%
\begin{align}
\ell\left(  \sigma\right)   &  =\left(  \text{the number of inversions of
}\sigma\right)  \ \ \ \ \ \ \ \ \ \ \left(  \text{by the definition of }%
\ell\left(  \sigma\right)  \right) \nonumber\\
&  =\left(  \text{the number of elements of }\operatorname*{Inv}\left(
\sigma\right)  \right) \nonumber\\
&  \ \ \ \ \ \ \ \ \ \ \left(  \text{since }\operatorname*{Inv}\left(
\sigma\right)  \text{ is the set of all inversions of }\sigma\right)
\nonumber\\
&  =\left\vert \operatorname*{Inv}\left(  \sigma\right)  \right\vert .
\label{sol.ps2.2.5.f.1}%
\end{align}

\textbf{(d)} Let $\sigma\in S_{n}$ be a permutation satisfying $\sigma\left(
1\right)  \leq\sigma\left(  2\right)  \leq\cdots\leq\sigma\left(  n\right)  $.
Then, Proposition \ref{prop.sol.ps2.2.5.d} shows that $\sigma
=\operatorname*{id}$. This solves Exercise \ref{exe.ps2.2.5} \textbf{(d)}.

\textbf{(f)} Let $\sigma\in S_{n}$. For every $\left(  i,j\right)
\in\operatorname*{Inv}\left(  \sigma\right)  $, we have $\left(  \sigma\left(
j\right)  ,\sigma\left(  i\right)  \right)  \in\operatorname*{Inv}\left(
\sigma^{-1}\right)  $\ \ \ \ \footnote{\textit{Proof.} Let $\left(
i,j\right)  \in\operatorname*{Inv}\left(  \sigma\right)  $. Then, $\left(
i,j\right)  $ is an inversion of $\sigma$ (since $\operatorname*{Inv}\left(
\sigma\right)  $ is the set of all inversions of $\sigma$). In other words,
$\left(  i,j\right)  $ is a pair of integers satisfying $1\leq i<j\leq n$ and
$\sigma\left(  i\right)  >\sigma\left(  j\right)  $ (by the definition of an
\textquotedblleft inversion of $\sigma$\textquotedblright). Hence,
$\sigma\left(  j\right)  <\sigma\left(  i\right)  $, so that $1\leq
\sigma\left(  j\right)  <\sigma\left(  i\right)  \leq n$; also, $\sigma
^{-1}\left(  \sigma\left(  i\right)  \right)  =i<j=\sigma^{-1}\left(
\sigma\left(  j\right)  \right)  $, so that $\sigma^{-1}\left(  \sigma\left(
j\right)  \right)  >\sigma^{-1}\left(  \sigma\left(  i\right)  \right)  $.
Therefore, $\left(  \sigma\left(  j\right)  ,\sigma\left(  i\right)  \right)
$ is a pair of integers $\left(  u,v\right)  $ satisfying $1\leq u<v\leq n$
and $\sigma^{-1}\left(  u\right)  >\sigma^{-1}\left(  v\right)  $ (since
$1\leq\sigma\left(  j\right)  <\sigma\left(  i\right)  \leq n$ and
$\sigma^{-1}\left(  \sigma\left(  j\right)  \right)  >\sigma^{-1}\left(
\sigma\left(  i\right)  \right)  $). In other words, $\left(  \sigma\left(
j\right)  ,\sigma\left(  i\right)  \right)  $ is an inversion of $\sigma^{-1}$
(because inversions of $\sigma^{-1}$ are defined as pairs of integers $\left(
u,v\right)  $ satisfying $1\leq u<v\leq n$ and $\sigma^{-1}\left(  u\right)
>\sigma^{-1}\left(  v\right)  $). In other words, $\left(  \sigma\left(
j\right)  ,\sigma\left(  i\right)  \right)  \in\operatorname*{Inv}\left(
\sigma^{-1}\right)  $ (since $\operatorname*{Inv}\left(  \sigma^{-1}\right)  $
is the set of all inversions of $\sigma^{-1}$), qed.}. Hence, we can define a
map%
\begin{align*}
\Phi:\operatorname*{Inv}\left(  \sigma\right)   &  \rightarrow
\operatorname*{Inv}\left(  \sigma^{-1}\right)  ,\\
\left(  i,j\right)   &  \mapsto\left(  \sigma\left(  j\right)  ,\sigma\left(
i\right)  \right)  .
\end{align*}
This map $\Phi$ is injective\footnote{\textit{Proof.} We simply need to prove
that an element $\left(  i,j\right)  \in\operatorname*{Inv}\left(
\sigma\right)  $ can be reconstructed from its image $\left(  \sigma\left(
j\right)  ,\sigma\left(  i\right)  \right)  $. But this is easy: If you know
$\left(  \sigma\left(  j\right)  ,\sigma\left(  i\right)  \right)  $, then you
know $\sigma\left(  j\right)  $ and $\sigma\left(  i\right)  $, and therefore
also $i$ (since $i=\sigma^{-1}\left(  \sigma\left(  i\right)  \right)  $) and
$j$ (since $j=\sigma^{-1}\left(  \sigma\left(  j\right)  \right)  $), and thus
also $\left(  i,j\right)  $.}. Thus, we have found an injective map from
$\operatorname*{Inv}\left(  \sigma\right)  $ to $\operatorname*{Inv}\left(
\sigma^{-1}\right)  $. Conversely, $\left\vert \operatorname*{Inv}\left(
\sigma\right)  \right\vert \leq\left\vert \operatorname*{Inv}\left(
\sigma^{-1}\right)  \right\vert $. But $\ell\left(  \sigma^{-1}\right)
=\left\vert \operatorname*{Inv}\left(  \sigma^{-1}\right)  \right\vert $ (by
(\ref{sol.ps2.2.5.f.1}), applied to $\sigma^{-1}$ instead of $\sigma$). Now,
(\ref{sol.ps2.2.5.f.1}) yields $\ell\left(  \sigma\right)  =\left\vert
\operatorname*{Inv}\left(  \sigma\right)  \right\vert \leq\left\vert
\operatorname*{Inv}\left(  \sigma^{-1}\right)  \right\vert =\ell\left(
\sigma^{-1}\right)  $.

Now, let us forget that we fixed $\sigma$. We thus have proven that
\begin{equation}
\ell\left(  \sigma\right)  \leq\ell\left(  \sigma^{-1}\right)
\ \ \ \ \ \ \ \ \ \ \text{for every }\sigma\in S_{n}. \label{sol.ps2.2.5.f.6}%
\end{equation}

Now, let $\sigma\in S_{n}$ again. We can apply (\ref{sol.ps2.2.5.f.6}) to
$\sigma^{-1}$ instead of $\sigma$, and thus obtain $\ell\left(  \sigma
^{-1}\right)  \leq\ell\left(  \underbrace{\left(  \sigma^{-1}\right)  ^{-1}%
}_{=\sigma}\right)  =\ell\left(  \sigma\right)  $. Combined with
(\ref{sol.ps2.2.5.f.6}), this yields $\ell\left(  \sigma\right)  =\ell\left(
\sigma^{-1}\right)  $. This solves Exercise \ref{exe.ps2.2.5} \textbf{(f)}.

\textbf{(a)} As I warned above, this solution will be a tedious formalization
of the argument sketched in Example \ref{exa.2.5}.

Let us first show a very simple fact: If $u$ and $v$ are two integers such
that $1\leq u<v\leq n$, and if $k\in\left\{  1,2,\ldots,n-1\right\}  $ is such
that $\left(  u,v\right)  \neq\left(  k,k+1\right)  $, then%
\begin{equation}
s_{k}\left(  u\right)  <s_{k}\left(  v\right)  \label{sol.ps2.2.5.a.sk-inc}%
\end{equation}
\footnote{\textit{Proof of (\ref{sol.ps2.2.5.a.sk-inc}):} We can prove
(\ref{sol.ps2.2.5.a.sk-inc}) by analyzing three cases (Case 1 is when $u=k$,
Case 2 is when $u=k+1$, and Case 3 is when $u\notin\left\{  k,k+1\right\}  $),
each of which can be split into three subcases (Subcase 1 is when $v=k$,
Subcase 2 is when $v=k+1$, and Subcase 3 is when $v\notin\left\{
k,k+1\right\}  $). These are (altogether) nine subcases, but four of them
(namely, Subcases 1 and 2 in Case 1, and Subcases 1 and 2 in Case 2) are
impossible (because $u<v$ and $\left(  u,v\right)  \neq\left(  k,k+1\right)
$), and the proof of (\ref{sol.ps2.2.5.a.sk-inc}) is easy in the remaining
five subcases.
\par
Here is a smarter way to prove (\ref{sol.ps2.2.5.a.sk-inc}): Let $u$ and $v$
be two integers such that $1\leq u<v\leq n$. Let $k\in\left\{  1,2,\ldots
,n-1\right\}  $ be such that $\left(  u,v\right)  \neq\left(  k,k+1\right)  $.
We need to prove (\ref{sol.ps2.2.5.a.sk-inc}). Indeed, assume the contrary.
Thus, $s_{k}\left(  u\right)  \geq s_{k}\left(  v\right)  $.
\par
But $u<v$ and thus $u\neq v$. The map $s_{k}$ is a permutation, thus
bijective, and therefore injective. Hence, $s_{k}\left(  u\right)  \neq
s_{k}\left(  v\right)  $ (since $u\neq v$). Combined with $s_{k}\left(
u\right)  \geq s_{k}\left(  v\right)  $, this yields $s_{k}\left(  u\right)
>s_{k}\left(  v\right)  $. Thus, $s_{k}\left(  u\right)  \geq s_{k}\left(
v\right)  +1$ (since $s_{k}\left(  u\right)  $ and $s_{k}\left(  v\right)  $
are integers), so that $s_{k}\left(  v\right)  +1\leq s_{k}\left(  u\right)
$.
\par
We have $u<v$ and thus $u+1\leq v$ (since $u$ and $v$ are integers).
\par
Recall that $s_{k}$ is the permutation in $S_{n}$ which swaps $k$ and $k+1$,
while leaving all other elements of $\left\{  1,2,\ldots,n\right\}  $
unchanged. Hence,
\begin{equation}
s_{k}\left(  p\right)  \leq p+1\ \ \ \ \ \ \ \ \ \ \text{for every }%
p\in\left\{  1,2,\ldots,n\right\}  , \label{sol.ps2.2.5.a.sk-inc.pf.1}%
\end{equation}
and this inequality becomes an equality only for $p=k$. For the same reason,
we have%
\begin{equation}
s_{k}\left(  p\right)  \geq p-1\ \ \ \ \ \ \ \ \ \ \text{for every }%
p\in\left\{  1,2,\ldots,n\right\}  , \label{sol.ps2.2.5.a.sk-inc.pf.2}%
\end{equation}
and this inequality becomes an equality only for $p=k+1$.
\par
Applying (\ref{sol.ps2.2.5.a.sk-inc.pf.1}) to $p=u$, we obtain $s_{k}\left(
u\right)  \leq u+1$. Applying (\ref{sol.ps2.2.5.a.sk-inc.pf.2}) to $p=v$, we
obtain $s_{k}\left(  v\right)  \geq v-1$, so that $v-1\leq s_{k}\left(
v\right)  $ and thus $v\leq s_{k}\left(  v\right)  +1\leq s_{k}\left(
u\right)  \leq u+1\leq v$.
\par
Combining $v\leq s_{k}\left(  v\right)  +1$ with $s_{k}\left(  v\right)
+1\leq v$, we obtain $v=s_{k}\left(  v\right)  +1$, so that $s_{k}\left(
v\right)  =v-1$. In other words, $s_{k}\left(  p\right)  =p-1$ holds for
$p=v$. But recall that the inequality (\ref{sol.ps2.2.5.a.sk-inc.pf.2})
becomes an equality only for $p=k+1$. In other words, $s_{k}\left(  p\right)
=p-1$ holds only for $p=k+1$. Applying this to $p=v$, we obtain $v=k+1$ (since
$s_{k}\left(  p\right)  =p-1$ holds for $p=v$).
\par
Combining $s_{k}\left(  u\right)  \leq u+1$ with $u+1\leq v\leq s_{k}\left(
u\right)  $, we obtain $s_{k}\left(  u\right)  =u+1$. In other words,
$s_{k}\left(  p\right)  =p+1$ holds for $p=u$. But recall that the inequality
(\ref{sol.ps2.2.5.a.sk-inc.pf.1}) becomes an equality only for $p=k$. In other
words, $s_{k}\left(  p\right)  =p+1$ holds only for $p=k$. Applying this to
$p=u$, we obtain $u=k$ (since $s_{k}\left(  p\right)  =p+1$ holds for $p=u$).
\par
Now, $\left(  \underbrace{u}_{=k},\underbrace{v}_{=k+1}\right)  =\left(
k,k+1\right)  $ contradicts $\left(  u,v\right)  \neq\left(  k,k+1\right)  $.
This contradiction proves that our assumption was wrong. Hence,
(\ref{sol.ps2.2.5.a.sk-inc}) is proven.}.

Recall that $s_{i}^{2}=\operatorname*{id}$ for each $i\in\left\{
1,2,\ldots,n-1\right\}  $. Applying this to $i=k$, we obtain $s_{k}%
^{2}=\operatorname*{id}$; thus, $s_{k}\circ s_{k}=s_{k}^{2}=\operatorname*{id}%
$ and therefore $s_{k}^{-1}=s_{k}$.

We shall now show that%
\begin{align}
&  \operatorname*{Inv}\left(  s_{k}\circ\sigma\right)  \setminus\left\{
\left(  \sigma^{-1}\left(  k\right)  ,\sigma^{-1}\left(  k+1\right)  \right)
\right\} \nonumber\\
&  =\operatorname*{Inv}\left(  \sigma\right)  \setminus\left\{  \left(
\sigma^{-1}\left(  k+1\right)  ,\sigma^{-1}\left(  k\right)  \right)
\right\}  \label{sol.ps2.2.5.a.equalsets}%
\end{align}
for every $\sigma\in S_{n}$ and $k\in\left\{  1,2,\ldots,n-1\right\}  $.

[Notice that we do not necessarily have $\left(  \sigma^{-1}\left(
k+1\right)  ,\sigma^{-1}\left(  k\right)  \right)  \in\operatorname*{Inv}%
\left(  \sigma\right)  $; nor do we always have $\left(  \sigma^{-1}\left(
k\right)  ,\sigma^{-1}\left(  k+1\right)  \right)  \in\operatorname*{Inv}%
\left(  s_{k}\circ\sigma\right)  $. In fact, for each given $\sigma$ and $k$,
exactly one of these two statements holds. But we can form the difference
$A\setminus B$ of two sets $A$ and $B$ even if $B$ is not a subset of $A$, so
the statement (\ref{sol.ps2.2.5.a.equalsets}) still makes sense.]

[\textit{Proof of (\ref{sol.ps2.2.5.a.equalsets}):} Let $\sigma\in S_{n}$ and
$k\in\left\{  1,2,\ldots,n-1\right\}  $. Let \newline$\left(  i,j\right)
\in\operatorname*{Inv}\left(  \sigma\right)  \setminus\left\{  \left(
\sigma^{-1}\left(  k+1\right)  ,\sigma^{-1}\left(  k\right)  \right)
\right\}  $. Thus, $\left(  i,j\right)  \in\operatorname*{Inv}\left(
\sigma\right)  $ and \newline$\left(  i,j\right)  \neq\left(  \sigma
^{-1}\left(  k+1\right)  ,\sigma^{-1}\left(  k\right)  \right)  $. Therefore,
$\left(  \sigma\left(  j\right)  ,\sigma\left(  i\right)  \right)  \neq\left(
k,k+1\right)  $\ \ \ \ \footnote{\textit{Proof.} Assume the contrary. Thus,
$\left(  \sigma\left(  j\right)  ,\sigma\left(  i\right)  \right)  =\left(
k,k+1\right)  $. Hence, $\sigma\left(  j\right)  =k$ and $\sigma\left(
i\right)  =k+1$. Hence, $\left(  \underbrace{\sigma^{-1}\left(  k+1\right)
}_{\substack{=i\\\text{(since }\sigma\left(  i\right)  =k+1\text{)}%
}},\underbrace{\sigma^{-1}\left(  k\right)  }_{\substack{=j\\\text{(since
}\sigma\left(  j\right)  =k\text{)}}}\right)  =\left(  i,j\right)  \neq\left(
\sigma^{-1}\left(  k+1\right)  ,\sigma^{-1}\left(  k\right)  \right)  $, which
is absurd. Hence, we have found a contradiction, so that our assumption was
wrong, qed.}.

We have $\left(  i,j\right)  \in\operatorname*{Inv}\left(  \sigma\right)  $.
In other words, $\left(  i,j\right)  $ is an inversion of $\sigma$. In other
words, $\left(  i,j\right)  $ is a pair of integers satisfying $1\leq i<j\leq
n$ and $\sigma\left(  i\right)  >\sigma\left(  j\right)  $. Now,
$\sigma\left(  j\right)  <\sigma\left(  i\right)  $, so that $1\leq
\sigma\left(  j\right)  <\sigma\left(  i\right)  \leq n$. Also, as we know,
$\left(  \sigma\left(  j\right)  ,\sigma\left(  i\right)  \right)  \neq\left(
k,k+1\right)  $. Hence, $s_{k}\left(  \sigma\left(  j\right)  \right)
<s_{k}\left(  \sigma\left(  i\right)  \right)  $ (by
(\ref{sol.ps2.2.5.a.sk-inc}), applied to $u=\sigma\left(  j\right)  $ and
$v=\sigma\left(  i\right)  $). Thus, $\left(  s_{k}\circ\sigma\right)  \left(
j\right)  =s_{k}\left(  \sigma\left(  j\right)  \right)  <s_{k}\left(
\sigma\left(  i\right)  \right)  =\left(  s_{k}\circ\sigma\right)  \left(
i\right)  $, hence $\left(  s_{k}\circ\sigma\right)  \left(  i\right)
>\left(  s_{k}\circ\sigma\right)  \left(  j\right)  $. Hence, $\left(
i,j\right)  $ is a pair of integers satisfying $1\leq i<j\leq n$ and $\left(
s_{k}\circ\sigma\right)  \left(  i\right)  >\left(  s_{k}\circ\sigma\right)
\left(  j\right)  $. In other words, $\left(  i,j\right)  $ is an inversion of
$s_{k}\circ\sigma$ (by the definition of an inversion). In other words,
$\left(  i,j\right)  \in\operatorname*{Inv}\left(  s_{k}\circ\sigma\right)  $
(since $\operatorname*{Inv}\left(  s_{k}\circ\sigma\right)  $ is defined as
the set of all inversions of $s_{k}\circ\sigma$). Furthermore, $\left(
i,j\right)  \neq\left(  \sigma^{-1}\left(  k\right)  ,\sigma^{-1}\left(
k+1\right)  \right)  $\ \ \ \ \footnote{\textit{Proof.} Assume the contrary.
Thus, $\left(  i,j\right)  =\left(  \sigma^{-1}\left(  k\right)  ,\sigma
^{-1}\left(  k+1\right)  \right)  $. Thus, $i=\sigma^{-1}\left(  k\right)  $
and $j=\sigma^{-1}\left(  k+1\right)  $. Hence, $\sigma\left(  i\right)  =k$
(since $i=\sigma^{-1}\left(  k\right)  $) and $\sigma\left(  j\right)  =k+1$
(since $j=\sigma^{-1}\left(  k+1\right)  $). Hence, $k+1=\sigma\left(
j\right)  <\sigma\left(  i\right)  =k<k+1$, which is absurd. Thus, we have
found a contradiction, so that our assumption must have been wrong, qed.}.
Thus, $\left(  i,j\right)  \in\operatorname*{Inv}\left(  s_{k}\circ
\sigma\right)  \setminus\left\{  \left(  \sigma^{-1}\left(  k\right)
,\sigma^{-1}\left(  k+1\right)  \right)  \right\}  $ (since $\left(
i,j\right)  \in\operatorname*{Inv}\left(  s_{k}\circ\sigma\right)  $ and
$\left(  i,j\right)  \neq\left(  \sigma^{-1}\left(  k\right)  ,\sigma
^{-1}\left(  k+1\right)  \right)  $).

Now, let us forget that we fixed $\left(  i,j\right)  $. We thus have shown
that every $\left(  i,j\right)  \in\operatorname*{Inv}\left(  \sigma\right)
\setminus\left\{  \left(  \sigma^{-1}\left(  k+1\right)  ,\sigma^{-1}\left(
k\right)  \right)  \right\}  $ satisfies \newline$\left(  i,j\right)
\in\operatorname*{Inv}\left(  s_{k}\circ\sigma\right)  \setminus\left\{
\left(  \sigma^{-1}\left(  k\right)  ,\sigma^{-1}\left(  k+1\right)  \right)
\right\}  $. In other words,%
\begin{align}
&  \operatorname*{Inv}\left(  \sigma\right)  \setminus\left\{  \left(
\sigma^{-1}\left(  k+1\right)  ,\sigma^{-1}\left(  k\right)  \right)  \right\}
\nonumber\\
&  \subseteq\operatorname*{Inv}\left(  s_{k}\circ\sigma\right)  \setminus
\left\{  \left(  \sigma^{-1}\left(  k\right)  ,\sigma^{-1}\left(  k+1\right)
\right)  \right\}  . \label{sol.ps2.2.5.a.equalsets.pf.1}%
\end{align}

Now, let $\tau=s_{k}\circ\sigma$. Then, $s_{k}\circ\underbrace{\tau}%
_{=s_{k}\circ\sigma}=\underbrace{s_{k}\circ s_{k}}_{=s_{k}^{2}%
=\operatorname*{id}}\circ\sigma=\operatorname*{id}\circ\sigma=\sigma$.
Moreover, $\tau^{-1}\left(  k\right)  =\sigma^{-1}\left(  k+1\right)
$\ \ \ \ \footnote{\textit{Proof.} We have $\underbrace{\tau}_{=s_{k}%
\circ\sigma}\left(  \sigma^{-1}\left(  k+1\right)  \right)  =\left(
s_{k}\circ\sigma\right)  \left(  \sigma^{-1}\left(  k+1\right)  \right)
=s_{k}\left(  \underbrace{\sigma\left(  \sigma^{-1}\left(  k+1\right)
\right)  }_{=k+1}\right)  =s_{k}\left(  k+1\right)  =k$ (by the definition of
$s_{k}$). Thus, $\tau^{-1}\left(  k\right)  =\sigma^{-1}\left(  k+1\right)  $,
qed.} and $\tau^{-1}\left(  k+1\right)  =\sigma^{-1}\left(  k\right)
$\ \ \ \ \footnote{for similar reasons}.

But recall that we have proven (\ref{sol.ps2.2.5.a.equalsets.pf.1}). The same
arguments, but carried out for $\tau$ instead of $\sigma$, show that%
\[
\operatorname*{Inv}\left(  \tau\right)  \setminus\left\{  \left(  \tau
^{-1}\left(  k+1\right)  ,\tau^{-1}\left(  k\right)  \right)  \right\}
\subseteq\operatorname*{Inv}\left(  s_{k}\circ\tau\right)  \setminus\left\{
\left(  \tau^{-1}\left(  k\right)  ,\tau^{-1}\left(  k+1\right)  \right)
\right\}  .
\]
Using the identities $s_{k}\circ\tau=\sigma$ and $\tau^{-1}\left(  k\right)
=\sigma^{-1}\left(  k+1\right)  $ and $\tau^{-1}\left(  k+1\right)
=\sigma^{-1}\left(  k\right)  $, we can rewrite this as follows:%
\[
\operatorname*{Inv}\left(  \tau\right)  \setminus\left\{  \left(  \sigma
^{-1}\left(  k\right)  ,\sigma^{-1}\left(  k+1\right)  \right)  \right\}
\subseteq\operatorname*{Inv}\left(  \sigma\right)  \setminus\left\{  \left(
\sigma^{-1}\left(  k+1\right)  ,\sigma^{-1}\left(  k\right)  \right)
\right\}  .
\]
Since $\tau=s_{k}\circ\sigma$, this further rewrites as follows:%
\[
\operatorname*{Inv}\left(  s_{k}\circ\sigma\right)  \setminus\left\{  \left(
\sigma^{-1}\left(  k\right)  ,\sigma^{-1}\left(  k+1\right)  \right)
\right\}  \subseteq\operatorname*{Inv}\left(  \sigma\right)  \setminus\left\{
\left(  \sigma^{-1}\left(  k+1\right)  ,\sigma^{-1}\left(  k\right)  \right)
\right\}  .
\]
Combining this with (\ref{sol.ps2.2.5.a.equalsets.pf.1}), we obtain%
\[
\operatorname*{Inv}\left(  s_{k}\circ\sigma\right)  \setminus\left\{  \left(
\sigma^{-1}\left(  k\right)  ,\sigma^{-1}\left(  k+1\right)  \right)
\right\}  =\operatorname*{Inv}\left(  \sigma\right)  \setminus\left\{  \left(
\sigma^{-1}\left(  k+1\right)  ,\sigma^{-1}\left(  k\right)  \right)
\right\}  .
\]
This proves (\ref{sol.ps2.2.5.a.equalsets}).]

Now, let $k\in\left\{  1,2,\ldots,n-1\right\}  $.

We shall first show that for every $\sigma\in S_{n}$, we have%
\begin{equation}
\ell\left(  s_{k}\circ\sigma\right)  =\ell\left(  \sigma\right)
+1\ \ \ \ \ \ \ \ \ \ \text{if }\sigma^{-1}\left(  k\right)  <\sigma
^{-1}\left(  k+1\right)  . \label{sol.ps2.2.5.a.part1}%
\end{equation}

[\textit{Proof of (\ref{sol.ps2.2.5.a.part1}):} Let $\sigma\in S_{n}$. Assume
that $\sigma^{-1}\left(  k\right)  <\sigma^{-1}\left(  k+1\right)  $. Then,
\newline$\left(  \sigma^{-1}\left(  k\right)  ,\sigma^{-1}\left(  k+1\right)
\right)  $ is a pair of integers satisfying $1\leq\sigma^{-1}\left(  k\right)
<\sigma^{-1}\left(  k+1\right)  \leq n$ and $\left(  s_{k}\circ\sigma\right)
\left(  \sigma^{-1}\left(  k\right)  \right)  >\left(  s_{k}\circ
\sigma\right)  \left(  \sigma^{-1}\left(  k+1\right)  \right)  $%
\ \ \ \ \footnote{because $\left(  s_{k}\circ\sigma\right)  \left(
\sigma^{-1}\left(  k\right)  \right)  =s_{k}\left(  \underbrace{\sigma\left(
\sigma^{-1}\left(  k\right)  \right)  }_{=k}\right)  =s_{k}\left(  k\right)
=k+1$ and similarly $\left(  s_{k}\circ\sigma\right)  \left(  \sigma
^{-1}\left(  k+1\right)  \right)  =k$, so that $\left(  s_{k}\circ
\sigma\right)  \left(  \sigma^{-1}\left(  k\right)  \right)  =k+1>k=\left(
s_{k}\circ\sigma\right)  \left(  \sigma^{-1}\left(  k+1\right)  \right)  $}.
In other words, \newline$\left(  \sigma^{-1}\left(  k\right)  ,\sigma
^{-1}\left(  k+1\right)  \right)  $ is an inversion of $s_{k}\circ\sigma$. In
other words, $\left(  \sigma^{-1}\left(  k\right)  ,\sigma^{-1}\left(
k+1\right)  \right)  \in\operatorname*{Inv}\left(  s_{k}\circ\sigma\right)  $.
Hence,%
\begin{equation}
\left\vert \operatorname*{Inv}\left(  s_{k}\circ\sigma\right)  \setminus
\left\{  \left(  \sigma^{-1}\left(  k\right)  ,\sigma^{-1}\left(  k+1\right)
\right)  \right\}  \right\vert =\left\vert \operatorname*{Inv}\left(
s_{k}\circ\sigma\right)  \right\vert -1 \label{sol.ps2.2.5.a.part1.pf.1}%
\end{equation}

On the other hand, $\left(  \sigma^{-1}\left(  k+1\right)  ,\sigma^{-1}\left(
k\right)  \right)  $ is not an inversion of $\sigma$ (because if it was an
inversion of $\sigma$, then we would have $1\leq\sigma^{-1}\left(  k+1\right)
<\sigma^{-1}\left(  k\right)  \leq n$ and therefore $\sigma^{-1}\left(
k+1\right)  <\sigma^{-1}\left(  k\right)  <\sigma^{-1}\left(  k+1\right)  $,
which would be absurd). In other words, $\left(  \sigma^{-1}\left(
k+1\right)  ,\sigma^{-1}\left(  k\right)  \right)  \notin\operatorname*{Inv}%
\left(  \sigma\right)  $. Thus,%
\[
\operatorname*{Inv}\left(  \sigma\right)  \setminus\left\{  \left(
\sigma^{-1}\left(  k+1\right)  ,\sigma^{-1}\left(  k\right)  \right)
\right\}  =\operatorname*{Inv}\left(  \sigma\right)  ,
\]
so that%
\begin{align}
\operatorname*{Inv}\left(  \sigma\right)   &  =\operatorname*{Inv}\left(
\sigma\right)  \setminus\left\{  \left(  \sigma^{-1}\left(  k+1\right)
,\sigma^{-1}\left(  k\right)  \right)  \right\} \nonumber\\
&  =\operatorname*{Inv}\left(  s_{k}\circ\sigma\right)  \setminus\left\{
\left(  \sigma^{-1}\left(  k\right)  ,\sigma^{-1}\left(  k+1\right)  \right)
\right\}  \ \ \ \ \ \ \ \ \ \ \left(  \text{by (\ref{sol.ps2.2.5.a.equalsets}%
)}\right)  . \label{sol.ps2.2.5.a.part1.pf.4}%
\end{align}
Now, (\ref{sol.ps2.2.5.f.1}) yields%
\begin{align}
\ell\left(  \sigma\right)   &  =\left\vert \operatorname*{Inv}\left(
\sigma\right)  \right\vert =\left\vert \operatorname*{Inv}\left(  s_{k}%
\circ\sigma\right)  \setminus\left\{  \left(  \sigma^{-1}\left(  k\right)
,\sigma^{-1}\left(  k+1\right)  \right)  \right\}  \right\vert
\ \ \ \ \ \ \ \ \ \ \left(  \text{by (\ref{sol.ps2.2.5.a.part1.pf.4})}\right)
\nonumber\\
&  =\left\vert \operatorname*{Inv}\left(  s_{k}\circ\sigma\right)  \right\vert
-1\ \ \ \ \ \ \ \ \ \ \left(  \text{by (\ref{sol.ps2.2.5.a.part1.pf.1}%
)}\right)  . \label{sol.ps2.2.5.a.part1.pf.6}%
\end{align}
But (\ref{sol.ps2.2.5.f.1}) (applied to $s_{k}\circ\sigma$ instead of $\sigma
$) yields $\ell\left(  s_{k}\circ\sigma\right)  =\left\vert
\operatorname*{Inv}\left(  s_{k}\circ\sigma\right)  \right\vert $. Hence,
(\ref{sol.ps2.2.5.a.part1.pf.6}) becomes%
\[
\ell\left(  \sigma\right)  =\underbrace{\left\vert \operatorname*{Inv}\left(
s_{k}\circ\sigma\right)  \right\vert }_{=\ell\left(  s_{k}\circ\sigma\right)
}-1=\ell\left(  s_{k}\circ\sigma\right)  -1,
\]
so that $\ell\left(  s_{k}\circ\sigma\right)  =\ell\left(  \sigma\right)  +1$.
This proves (\ref{sol.ps2.2.5.a.part1}).]

Next, we will show that for every $\sigma\in S_{n}$, we have
\begin{equation}
\ell\left(  s_{k}\circ\sigma\right)  =\ell\left(  \sigma\right)
-1\ \ \ \ \ \ \ \ \ \ \text{if }\sigma^{-1}\left(  k\right)  >\sigma
^{-1}\left(  k+1\right)  . \label{sol.ps2.2.5.a.part2}%
\end{equation}

[\textit{Proof of (\ref{sol.ps2.2.5.a.part2}):} Let $\sigma\in S_{n}$. Assume
that $\sigma^{-1}\left(  k\right)  >\sigma^{-1}\left(  k+1\right)  $. But
$\sigma^{-1}\left(  k\right)  =\left(  s_{k}\circ\sigma\right)  ^{-1}\left(
k+1\right)  $\ \ \ \ \footnote{since $\left(  s_{k}\circ\sigma\right)  \left(
\sigma^{-1}\left(  k\right)  \right)  =s_{k}\left(  \underbrace{\sigma\left(
\sigma^{-1}\left(  k\right)  \right)  }_{=k}\right)  =s_{k}\left(  k\right)
=k+1$} and $\sigma^{-1}\left(  k+1\right)  =\left(  s_{k}\circ\sigma\right)
^{-1}\left(  k\right)  $\ \ \ \ \footnote{for similar reasons}. Thus, $\left(
s_{k}\circ\sigma\right)  ^{-1}\left(  k\right)  =\sigma^{-1}\left(
k+1\right)  <\sigma^{-1}\left(  k\right)  =\left(  s_{k}\circ\sigma\right)
^{-1}\left(  k+1\right)  $. Hence, we can apply (\ref{sol.ps2.2.5.a.part1}) to
$s_{k}\circ\sigma$ instead of $\sigma$. As a result, we obtain%
\[
\ell\left(  s_{k}\circ s_{k}\circ\sigma\right)  =\ell\left(  s_{k}\circ
\sigma\right)  +1.
\]
Since $\underbrace{s_{k}\circ s_{k}}_{=s_{k}^{2}=\operatorname*{id}}%
\circ\sigma=\operatorname*{id}\circ\sigma=\sigma$, this rewrites as
$\ell\left(  \sigma\right)  =\ell\left(  s_{k}\circ\sigma\right)  +1$, so that
$\ell\left(  s_{k}\circ\sigma\right)  =\ell\left(  \sigma\right)  -1$. This
proves (\ref{sol.ps2.2.5.a.part2}).]

\begin{vershort}
Now, (\ref{eq.exe.2.5.a.2}) follows immediately by combining
(\ref{sol.ps2.2.5.a.part1}) with (\ref{sol.ps2.2.5.a.part2}).\footnote{The
term \textquotedblleft$%
\begin{cases}
\ell\left(  \sigma\right)  +1, & \text{if }\sigma\left(  k\right)
<\sigma\left(  k+1\right)  ;\\
\ell\left(  \sigma\right)  -1, & \text{if }\sigma\left(  k\right)
>\sigma\left(  k+1\right)
\end{cases}
$\textquotedblright\ in (\ref{eq.exe.2.5.a.2}) makes sense because every
$\sigma\in S_{n}$ and every $k\in\left\{  1,2,\ldots,n-1\right\}  $ satisfies
exactly one of the conditions $\sigma\left(  k\right)  <\sigma\left(
k+1\right)  $ and $\sigma\left(  k\right)  >\sigma\left(  k+1\right)  $.
(Indeed, $\sigma\left(  k\right)  =\sigma\left(  k+1\right)  $ is impossible,
because every permutation $\sigma\in S_{n}$ is injective.)}
\end{vershort}

\begin{verlong}
We will soon prove (\ref{eq.exe.2.5.a.1}) and (\ref{eq.exe.2.5.a.2}). Let us
first show that the right-hand sides of (\ref{eq.exe.2.5.a.1}) and
(\ref{eq.exe.2.5.a.2}) are always well-defined. Indeed, for every $\sigma\in
S_{n}$ and every $k\in\left\{  1,2,\ldots,n-1\right\}  $, the right-hand side
of (\ref{eq.exe.2.5.a.1}) is well-defined\footnote{\textit{Proof.} Let
$\sigma\in S_{n}$ and $k\in\left\{  1,2,\ldots,n-1\right\}  $. The map
$\sigma$ is a permutation and thus injective. Hence, $\sigma\left(  k\right)
\neq\sigma\left(  k+1\right)  $ (since $k\neq k+1$). Thus, either
$\sigma\left(  k\right)  <\sigma\left(  k+1\right)  $ or $\sigma\left(
k\right)  >\sigma\left(  k+1\right)  $. More precisely, exactly one of the
conditions $\sigma\left(  k\right)  <\sigma\left(  k+1\right)  $ and
$\sigma\left(  k\right)  >\sigma\left(  k+1\right)  $ is satisfied. Hence, the
right-hand side of (\ref{eq.exe.2.5.a.1}) is well-defined.}, and the
right-hand side of (\ref{eq.exe.2.5.a.2}) is
well-defined\footnote{\textit{Proof.} Let $\sigma\in S_{n}$ and $k\in\left\{
1,2,\ldots,n-1\right\}  $. The map $\sigma^{-1}$ is a permutation and thus
injective. Hence, $\sigma^{-1}\left(  k\right)  \neq\sigma^{-1}\left(
k+1\right)  $ (since $k\neq k+1$). Thus, either $\sigma^{-1}\left(  k\right)
<\sigma^{-1}\left(  k+1\right)  $ or $\sigma^{-1}\left(  k\right)
>\sigma^{-1}\left(  k+1\right)  $. More precisely, exactly one of the
conditions $\sigma^{-1}\left(  k\right)  <\sigma^{-1}\left(  k+1\right)  $ and
$\sigma^{-1}\left(  k\right)  >\sigma^{-1}\left(  k+1\right)  $ is satisfied.
Hence, the right-hand side of (\ref{eq.exe.2.5.a.2}) is well-defined.}.

Now, for every $\sigma\in S_{n}$ and every $k\in\left\{  1,2,\ldots
,n-1\right\}  $, the equality (\ref{eq.exe.2.5.a.2}) holds (since it follows
immediately by combining (\ref{sol.ps2.2.5.a.part1}) with
(\ref{sol.ps2.2.5.a.part2})).
\end{verlong}

It remains to prove (\ref{eq.exe.2.5.a.1}). Indeed, let $\sigma\in S_{n}$. Let
us recall that $\left(  \alpha\circ\beta\right)  ^{-1}=\beta^{-1}\circ
\alpha^{-1}$ for any two permutations $\alpha$ and $\beta$ in $S_{n}$.
Applying this to $\alpha=s_{k}$ and $\beta=\sigma^{-1}$, we obtain $\left(
s_{k}\circ\sigma^{-1}\right)  ^{-1}=\underbrace{\left(  \sigma^{-1}\right)
^{-1}}_{=\sigma}\circ\underbrace{s_{k}^{-1}}_{=s_{k}}=\sigma\circ s_{k}$. But
Exercise \ref{exe.ps2.2.5} \textbf{(f)} yields $\ell\left(  \sigma\right)
=\ell\left(  \sigma^{-1}\right)  $. Also, Exercise \ref{exe.ps2.2.5}
\textbf{(f)} (applied to $s_{k}\circ\sigma^{-1}$ instead of $\sigma$) yields
$\ell\left(  s_{k}\circ\sigma^{-1}\right)  =\ell\left(  \underbrace{\left(
s_{k}\circ\sigma^{-1}\right)  ^{-1}}_{=\sigma\circ s_{k}}\right)  =\ell\left(
\sigma\circ s_{k}\right)  $. But applying (\ref{eq.exe.2.5.a.2}) to
$\sigma^{-1}$ instead of $\sigma$, we obtain%
\[
\ell\left(  s_{k}\circ\sigma^{-1}\right)  =
\begin{cases}
\ell\left(  \sigma^{-1}\right)  +1, & \text{if }\left(  \sigma^{-1}\right)
^{-1}\left(  k\right)  <\left(  \sigma^{-1}\right)  ^{-1}\left(  k+1\right)
;\\
\ell\left(  \sigma^{-1}\right)  -1, & \text{if }\left(  \sigma^{-1}\right)
^{-1}\left(  k\right)  >\left(  \sigma^{-1}\right)  ^{-1}\left(  k+1\right)
\end{cases}
.
\]
Since $\ell\left(  s_{k}\circ\sigma^{-1}\right)  =\ell\left(  \sigma\circ
s_{k}\right)  $, $\ell\left(  \sigma^{-1}\right)  =\ell\left(  \sigma\right)
$ and $\left(  \sigma^{-1}\right)  ^{-1}=\sigma$, this equality rewrites as
follows:%
\[
\ell\left(  \sigma\circ s_{k}\right)  =
\begin{cases}
\ell\left(  \sigma\right)  +1, & \text{if }\sigma\left(  k\right)
<\sigma\left(  k+1\right)  ;\\
\ell\left(  \sigma\right)  -1, & \text{if }\sigma\left(  k\right)
>\sigma\left(  k+1\right)
\end{cases}
.
\]
This proves (\ref{eq.exe.2.5.a.1}), and thus completes the solution of
Exercise \ref{exe.ps2.2.5} \textbf{(a)}.

\textbf{(e)} We shall solve Exercise \ref{exe.ps2.2.5} \textbf{(e)} by
induction over $\ell\left(  \sigma\right)  $:

\textit{Induction base:} Exercise \ref{exe.ps2.2.5} \textbf{(e)} holds in the
case when $\ell\left(  \sigma\right)  =0$\ \ \ \ \footnote{\textit{Proof.} Let
$\sigma\in S_{n}$ be such that $\ell\left(  \sigma\right)  =0$. We need to
show that $\sigma$ can be written as a composition of $\ell\left(
\sigma\right)  $ permutations of the form $s_{k}$ (with $k\in\left\{
1,2,\ldots,n-1\right\}  $).
\par
Recall that the composition of $0$ permutations in $S_{n}$ is
$\operatorname*{id}$ (by definition).
\par
We have $\ell\left(  \sigma\right)  =0$, and thus $\sigma=\operatorname*{id}$
(by Corollary \ref{cor.sol.ps2.2.5.d2}). Therefore, $\sigma$ is a composition
of $0$ permutations of the form $s_{k}$ (with $k\in\left\{  1,2,\ldots
,n-1\right\}  $) (because the composition of $0$ permutations in $S_{n}$ is
$\operatorname*{id}$). In other words, $\sigma$ is a composition of
$\ell\left(  \sigma\right)  $ permutations of the form $s_{k}$ (with
$k\in\left\{  1,2,\ldots,n-1\right\}  $) (since $\ell\left(  \sigma\right)
=0$). Thus, Exercise \ref{exe.ps2.2.5} \textbf{(e)} is solved in the case when
$\ell\left(  \sigma\right)  =0$, qed.}. This completes the induction base.

\textit{Induction step:} Let $L$ be a positive integer. Assume that Exercise
\ref{exe.ps2.2.5} \textbf{(e)} is solved in the case when $\ell\left(
\sigma\right)  =L-1$. We need to solve Exercise \ref{exe.ps2.2.5} \textbf{(e)}
in the case when $\ell\left(  \sigma\right)  =L$.

So let $\sigma\in S_{n}$ be such that $\ell\left(  \sigma\right)  =L$. We need
to show that $\sigma$ can be written as a composition of $\ell\left(
\sigma\right)  $ permutations of the form $s_{k}$ (with $k\in\left\{
1,2,\ldots,n-1\right\}  $).

There exists a $k\in\left\{  1,2,\ldots,n-1\right\}  $ such that
$\sigma\left(  k\right)  >\sigma\left(  k+1\right)  $%
\ \ \ \ \footnote{\textit{Proof.} Assume the contrary. Then, every
$k\in\left\{  1,2,\ldots,n-1\right\}  $ satisfies $\sigma\left(  k\right)
\leq\sigma\left(  k+1\right)  $. In other words, $\sigma\left(  1\right)
\leq\sigma\left(  2\right)  \leq\cdots\leq\sigma\left(  n\right)  $. Exercise
\ref{exe.ps2.2.5} \textbf{(d)} yields $\sigma=\operatorname*{id}$. Hence,
$\ell\left(  \sigma\right)  =\ell\left(  \operatorname*{id}\right)  =0$, so
that $0=\ell\left(  \sigma\right)  =L$. This contradicts the fact that $L$ is
a positive integer. This contradiction shows that our assumption was wrong,
qed.}. Fix such a $k$, and denote it by $j$. Thus, $j$ is an element of
$\left\{  1,2,\ldots,n-1\right\}  $ and satisfies $\sigma\left(  j\right)
>\sigma\left(  j+1\right)  $. From (\ref{eq.exe.2.5.a.1}) (applied to $k=j$),
we obtain%
\begin{align*}
\ell\left(  \sigma\circ s_{j}\right)   &  =
\begin{cases}
\ell\left(  \sigma\right)  +1, & \text{if }\sigma\left(  j\right)
<\sigma\left(  j+1\right)  ;\\
\ell\left(  \sigma\right)  -1, & \text{if }\sigma\left(  j\right)
>\sigma\left(  j+1\right)
\end{cases}
\\
&  =\underbrace{\ell\left(  \sigma\right)  }_{=L}-1\ \ \ \ \ \ \ \ \ \ \left(
\text{since }\sigma\left(  j\right)  >\sigma\left(  j+1\right)  \right) \\
&  =L-1.
\end{align*}
Hence, we can apply Exercise \ref{exe.ps2.2.5} \textbf{(e)} to $\sigma\circ
s_{j}$ instead of $\sigma$ (because we assumed that Exercise \ref{exe.ps2.2.5}
\textbf{(e)} is solved in the case when $\ell\left(  \sigma\right)  =L-1$). As
a result, we conclude that $\sigma\circ s_{j}$ can be written as a composition
of $\ell\left(  \sigma\circ s_{j}\right)  $ permutations of the form $s_{k}$
(with $k\in\left\{  1,2,\ldots,n-1\right\}  $). In other words, $\sigma\circ
s_{j}$ can be written as a composition of $L-1$ permutations of the form
$s_{k}$ (with $k\in\left\{  1,2,\ldots,n-1\right\}  $) (since $\ell\left(
\sigma\circ s_{j}\right)  =L-1$). In other words, there exists an $\left(
L-1\right)  $-tuple $\left(  k_{1},k_{2},\ldots,k_{L-1}\right)  \in\left\{
1,2,\ldots,n-1\right\}  ^{L-1}$ such that $\sigma\circ s_{j}=s_{k_{1}}\circ
s_{k_{2}}\circ\cdots\circ s_{k_{L-1}}$. Consider this $\left(  k_{1}%
,k_{2},\ldots,k_{L-1}\right)  $.

We have $\sigma\circ\underbrace{s_{j}\circ s_{j}}_{=s_{j}^{2}%
=\operatorname*{id}}=\sigma$ and thus%
\[
\sigma=\underbrace{\sigma\circ s_{j}}_{=s_{k_{1}}\circ s_{k_{2}}\circ
\cdots\circ s_{k_{L-1}}}\circ s_{j}=s_{k_{1}}\circ s_{k_{2}}\circ\cdots\circ
s_{k_{L-1}}\circ s_{j}.
\]
The right hand side of this equality is a composition of $L$ permutations of
the form $s_{k}$ (with $k\in\left\{  1,2,\ldots,n-1\right\}  $). Thus,
$\sigma$ can be written as a composition of $L$ permutations of the form
$s_{k}$ (with $k\in\left\{  1,2,\ldots,n-1\right\}  $). In other words,
$\sigma$ can be written as a composition of $\ell\left(  \sigma\right)  $
permutations of the form $s_{k}$ (with $k\in\left\{  1,2,\ldots,n-1\right\}
$) (since $\ell\left(  \sigma\right)  =L$). This solves Exercise
\ref{exe.ps2.2.5} \textbf{(e)} in the case when $\ell\left(  \sigma\right)
=L$. The induction step is thus complete, and Exercise \ref{exe.ps2.2.5}
\textbf{(e)} is solved by induction.

\textbf{(b)} From (\ref{eq.exe.2.5.a.1}), we can easily conclude that%
\begin{equation}
\ell\left(  \sigma\circ s_{k}\right)  \equiv\ell\left(  \sigma\right)
+1\operatorname{mod}2 \label{sol.ps2.2.5.b.stepper}%
\end{equation}
for every $\sigma\in S_{n}$ and every $k\in\left\{  1,2,\ldots,n-1\right\}  $.

\begin{verlong}
[\textit{Proof of (\ref{sol.ps2.2.5.b.stepper}):} Let $\sigma\in S_{n}$ and
$k\in\left\{  1,2,\ldots,n-1\right\}  $. From (\ref{eq.exe.2.5.a.1}), we
obtain%
\begin{align*}
\ell\left(  \sigma\circ s_{k}\right)   &  =%
\begin{cases}
\ell\left(  \sigma\right)  +1, & \text{if }\sigma\left(  k\right)
<\sigma\left(  k+1\right)  ;\\
\ell\left(  \sigma\right)  -1, & \text{if }\sigma\left(  k\right)
>\sigma\left(  k+1\right)
\end{cases}
\\
&  \equiv%
\begin{cases}
\ell\left(  \sigma\right)  +1, & \text{if }\sigma\left(  k\right)
<\sigma\left(  k+1\right)  ;\\
\ell\left(  \sigma\right)  +1, & \text{if }\sigma\left(  k\right)
>\sigma\left(  k+1\right)
\end{cases}
\\
&  \ \ \ \ \ \ \ \ \ \ \left(  \text{since }\ell\left(  \sigma\right)
-1\equiv\ell\left(  \sigma\right)  +1\operatorname{mod}2\text{ in the case
when }\sigma\left(  k\right)  >\sigma\left(  k+1\right)  \right) \\
&  =\ell\left(  \sigma\right)  +1\operatorname{mod}2.
\end{align*}
This proves (\ref{sol.ps2.2.5.b.stepper}).]
\end{verlong}

Thus, using induction, it is easy to prove that%
\begin{equation}
\ell\left(  \sigma\circ\left(  s_{k_{1}}\circ s_{k_{2}}\circ\cdots\circ
s_{k_{p}}\right)  \right)  \equiv\ell\left(  \sigma\right)
+p\operatorname{mod}2 \label{sol.ps2.2.5.b.snake}%
\end{equation}
for every $\sigma\in S_{n}$, every $p\in\mathbb{N}$ and every $\left(
k_{1},k_{2},\ldots,k_{p}\right)  \in\left\{  1,2,\ldots,n-1\right\}  ^{p}$.

\begin{verlong}
[\textit{Proof of (\ref{sol.ps2.2.5.b.snake}):} Let $\sigma\in S_{n}$, let
$p\in\mathbb{N}$ and let $\left(  k_{1},k_{2},\ldots,k_{p}\right)  \in\left\{
1,2,\ldots,n-1\right\}  ^{p}$. We shall prove that%
\begin{equation}
\ell\left(  \sigma\circ\left(  s_{k_{1}}\circ s_{k_{2}}\circ\cdots\circ
s_{k_{q}}\right)  \right)  \equiv\ell\left(  \sigma\right)
+q\operatorname{mod}2 \label{sol.ps2.2.5.b.snake.pf.1}%
\end{equation}
for all $q\in\left\{  0,1,\ldots,p\right\}  $.

Indeed, let us prove (\ref{sol.ps2.2.5.b.snake.pf.1}) by induction over $q$.

\textit{Induction base:} We have $\ell\left(  \sigma\circ\underbrace{\left(
s_{k_{1}}\circ s_{k_{2}}\circ\cdots\circ s_{k_{0}}\right)  }_{=\left(
\text{empty composition}\right)  =\operatorname*{id}}\right)  =\ell\left(
\sigma\circ\operatorname*{id}\right)  =\ell\left(  \sigma\right)  =\ell\left(
\sigma\right)  +0\equiv\ell\left(  \sigma\right)  +0\operatorname{mod}2$. In
other words, (\ref{sol.ps2.2.5.b.snake.pf.1}) holds for $q=0$. This completes
the induction base.

\textit{Induction step:} Let $Q\in\left\{  0,1,\ldots,p\right\}  $ be
positive. Assume that (\ref{sol.ps2.2.5.b.snake.pf.1}) holds for $q=Q-1$. We
need to show that (\ref{sol.ps2.2.5.b.snake.pf.1}) holds for $q=Q$.

We have assumed that (\ref{sol.ps2.2.5.b.snake.pf.1}) holds for $q=Q-1$. In
other words, we have%
\[
\ell\left(  \sigma\circ\left(  s_{k_{1}}\circ s_{k_{2}}\circ\cdots\circ
s_{k_{Q-1}}\right)  \right)  \equiv\ell\left(  \sigma\right)  +\left(
Q-1\right)  \operatorname{mod}2.
\]
Now,%
\begin{align*}
&  \ell\left(  \sigma\circ\underbrace{\left(  s_{k_{1}}\circ s_{k_{2}}%
\circ\cdots\circ s_{k_{Q}}\right)  }_{=\left(  s_{k_{1}}\circ s_{k_{2}}%
\circ\cdots\circ s_{k_{Q-1}}\right)  \circ s_{k_{Q}}}\right) \\
&  =\ell\left(  \sigma\circ\left(  s_{k_{1}}\circ s_{k_{2}}\circ\cdots\circ
s_{k_{Q-1}}\right)  \circ s_{k_{Q}}\right) \\
&  \equiv\underbrace{\ell\left(  \sigma\circ\left(  s_{k_{1}}\circ s_{k_{2}%
}\circ\cdots\circ s_{k_{Q-1}}\right)  \right)  }_{\equiv\ell\left(
\sigma\right)  +\left(  Q-1\right)  \operatorname{mod}2}+1\\
&  \ \ \ \ \ \ \ \ \ \ \left(
\begin{array}
[c]{c}%
\text{by (\ref{sol.ps2.2.5.b.stepper}), applied to }\sigma\circ\left(
s_{k_{1}}\circ s_{k_{2}}\circ\cdots\circ s_{k_{Q-1}}\right)  \text{ and }%
k_{Q}\\
\text{instead of }\sigma\text{ and }k
\end{array}
\right) \\
&  \equiv\ell\left(  \sigma\right)  +\left(  Q-1\right)  +1=\ell\left(
\sigma\right)  +Q\operatorname{mod}2.
\end{align*}
In other words, (\ref{sol.ps2.2.5.b.snake.pf.1}) holds for $q=Q$. This
completes the induction step. Thus, (\ref{sol.ps2.2.5.b.snake.pf.1}) is proven
by induction.

Now, applying (\ref{sol.ps2.2.5.b.snake.pf.1}) to $q=p$, we obtain
$\ell\left(  \sigma\circ\left(  s_{k_{1}}\circ s_{k_{2}}\circ\cdots\circ
s_{k_{p}}\right)  \right)  \equiv\ell\left(  \sigma\right)
+p\operatorname{mod}2$. This proves (\ref{sol.ps2.2.5.b.snake}).]
\end{verlong}

Now, let $\sigma$ and $\tau$ be two permutations in $S_{n}$. Exercise
\ref{exe.ps2.2.5} \textbf{(e)} (applied to $\tau$ instead of $\sigma$) yields
that $\tau$ can be written as a composition of $\ell\left(  \tau\right)  $
permutations of the form $s_{k}$ (with $k\in\left\{  1,2,\ldots,n-1\right\}
$). In other words, there exists an $\ell\left(  \tau\right)  $-tuple $\left(
k_{1},k_{2},\ldots,k_{\ell\left(  \tau\right)  }\right)  \in\left\{
1,2,\ldots,n-1\right\}  ^{\ell\left(  \tau\right)  }$ such that $\tau
=s_{k_{1}}\circ s_{k_{2}}\circ\cdots\circ s_{k_{\ell\left(  \tau\right)  }}$.
Consider this $\left(  k_{1},k_{2},\ldots,k_{\ell\left(  \tau\right)
}\right)  $. Then,%
\[
\ell\left(  \sigma\circ\underbrace{\tau}_{=s_{k_{1}}\circ s_{k_{2}}\circ
\cdots\circ s_{k_{\ell\left(  \tau\right)  }}}\right)  =\ell\left(
\sigma\circ\left(  s_{k_{1}}\circ s_{k_{2}}\circ\cdots\circ s_{k_{\ell\left(
\tau\right)  }}\right)  \right)  \equiv\ell\left(  \sigma\right)  +\ell\left(
\tau\right)  \operatorname{mod}2
\]
(by (\ref{sol.ps2.2.5.b.snake}), applied to $p=\ell\left(  \tau\right)  $).
This solves Exercise \ref{exe.ps2.2.5} \textbf{(b)}.

\textbf{(c)} The solution of Exercise \ref{exe.ps2.2.5} \textbf{(c)} is mostly
parallel to our above solution to Exercise \ref{exe.ps2.2.5} \textbf{(b)}.

From (\ref{eq.exe.2.5.a.1}), we can easily conclude that%
\begin{equation}
\ell\left(  \sigma\circ s_{k}\right)  \leq\ell\left(  \sigma\right)  +1
\label{sol.ps2.2.5.c.stepper}%
\end{equation}
for every $\sigma\in S_{n}$ and every $k\in\left\{  1,2,\ldots,n-1\right\}  $.

\begin{verlong}
[\textit{Proof of (\ref{sol.ps2.2.5.c.stepper}):} Let $\sigma\in S_{n}$ and
$k\in\left\{  1,2,\ldots,n-1\right\}  $. From (\ref{eq.exe.2.5.a.1}), we
obtain%
\begin{align*}
\ell\left(  \sigma\circ s_{k}\right)   &  =%
\begin{cases}
\ell\left(  \sigma\right)  +1, & \text{if }\sigma\left(  k\right)
<\sigma\left(  k+1\right)  ;\\
\ell\left(  \sigma\right)  -1, & \text{if }\sigma\left(  k\right)
>\sigma\left(  k+1\right)
\end{cases}
\\
&  \leq%
\begin{cases}
\ell\left(  \sigma\right)  +1, & \text{if }\sigma\left(  k\right)
<\sigma\left(  k+1\right)  ;\\
\ell\left(  \sigma\right)  +1, & \text{if }\sigma\left(  k\right)
>\sigma\left(  k+1\right)
\end{cases}
\\
&  \ \ \ \ \ \ \ \ \ \ \left(  \text{since }\ell\left(  \sigma\right)
-1\leq\ell\left(  \sigma\right)  +1\text{ in the case when }\sigma\left(
k\right)  >\sigma\left(  k+1\right)  \right) \\
&  =\ell\left(  \sigma\right)  +1.
\end{align*}
This proves (\ref{sol.ps2.2.5.c.stepper}).]
\end{verlong}

Thus, using induction, it is easy to prove that%
\begin{equation}
\ell\left(  \sigma\circ\left(  s_{k_{1}}\circ s_{k_{2}}\circ\cdots\circ
s_{k_{p}}\right)  \right)  \leq\ell\left(  \sigma\right)  +p
\label{sol.ps2.2.5.c.snake}%
\end{equation}
for every $\sigma\in S_{n}$, every $p\in\mathbb{N}$ and every $\left(
k_{1},k_{2},\ldots,k_{p}\right)  \in\left\{  1,2,\ldots,n-1\right\}  ^{p}$.

\begin{verlong}
[\textit{Proof of (\ref{sol.ps2.2.5.c.snake}):} Let $\sigma\in S_{n}$, let
$p\in\mathbb{N}$ and let $\left(  k_{1},k_{2},\ldots,k_{p}\right)  \in\left\{
1,2,\ldots,n-1\right\}  ^{p}$. We shall prove that%
\begin{equation}
\ell\left(  \sigma\circ\left(  s_{k_{1}}\circ s_{k_{2}}\circ\cdots\circ
s_{k_{q}}\right)  \right)  \leq\ell\left(  \sigma\right)  +q
\label{sol.ps2.2.5.c.snake.pf.1}%
\end{equation}
for all $q\in\left\{  0,1,\ldots,p\right\}  $.

Indeed, let us prove (\ref{sol.ps2.2.5.c.snake.pf.1}) by induction over $q$.

\textit{Induction base:} We have $\ell\left(  \sigma\circ\underbrace{\left(
s_{k_{1}}\circ s_{k_{2}}\circ\cdots\circ s_{k_{0}}\right)  }_{=\left(
\text{empty composition}\right)  =\operatorname*{id}}\right)  =\ell\left(
\sigma\circ\operatorname*{id}\right)  =\ell\left(  \sigma\right)  =\ell\left(
\sigma\right)  +0\leq\ell\left(  \sigma\right)  +0$. In other words,
(\ref{sol.ps2.2.5.c.snake.pf.1}) holds for $q=0$. This completes the induction base.

\textit{Induction step:} Let $Q\in\left\{  0,1,\ldots,p\right\}  $ be
positive. Assume that (\ref{sol.ps2.2.5.c.snake.pf.1}) holds for $q=Q-1$. We
need to show that (\ref{sol.ps2.2.5.c.snake.pf.1}) holds for $q=Q$.

We have assumed that (\ref{sol.ps2.2.5.c.snake.pf.1}) holds for $q=Q-1$. In
other words, we have%
\[
\ell\left(  \sigma\circ\left(  s_{k_{1}}\circ s_{k_{2}}\circ\cdots\circ
s_{k_{Q-1}}\right)  \right)  \leq\ell\left(  \sigma\right)  +\left(
Q-1\right)  .
\]
Now,%
\begin{align*}
&  \ell\left(  \sigma\circ\underbrace{\left(  s_{k_{1}}\circ s_{k_{2}}%
\circ\cdots\circ s_{k_{Q}}\right)  }_{=\left(  s_{k_{1}}\circ s_{k_{2}}%
\circ\cdots\circ s_{k_{Q-1}}\right)  \circ s_{k_{Q}}}\right) \\
&  =\ell\left(  \sigma\circ\left(  s_{k_{1}}\circ s_{k_{2}}\circ\cdots\circ
s_{k_{Q-1}}\right)  \circ s_{k_{Q}}\right) \\
&  \leq\underbrace{\ell\left(  \sigma\circ\left(  s_{k_{1}}\circ s_{k_{2}%
}\circ\cdots\circ s_{k_{Q-1}}\right)  \right)  }_{\leq\ell\left(
\sigma\right)  +\left(  Q-1\right)  }+1\\
&  \ \ \ \ \ \ \ \ \ \ \left(
\begin{array}
[c]{c}%
\text{by (\ref{sol.ps2.2.5.c.stepper}), applied to }\sigma\circ\left(
s_{k_{1}}\circ s_{k_{2}}\circ\cdots\circ s_{k_{Q-1}}\right)  \text{ and }%
k_{Q}\\
\text{instead of }\sigma\text{ and }k
\end{array}
\right) \\
&  \leq\ell\left(  \sigma\right)  +\left(  Q-1\right)  +1=\ell\left(
\sigma\right)  +Q.
\end{align*}
In other words, (\ref{sol.ps2.2.5.c.snake.pf.1}) holds for $q=Q$. This
completes the induction step. Thus, (\ref{sol.ps2.2.5.c.snake.pf.1}) is proven
by induction.

Now, applying (\ref{sol.ps2.2.5.c.snake.pf.1}) to $q=p$, we obtain
$\ell\left(  \sigma\circ\left(  s_{k_{1}}\circ s_{k_{2}}\circ\cdots\circ
s_{k_{p}}\right)  \right)  \leq\ell\left(  \sigma\right)  +p$. This proves
(\ref{sol.ps2.2.5.c.snake}).]
\end{verlong}

Now, let $\sigma$ and $\tau$ be two permutations in $S_{n}$. Exercise
\ref{exe.ps2.2.5} \textbf{(e)} (applied to $\tau$ instead of $\sigma$) yields
that $\tau$ can be written as a composition of $\ell\left(  \tau\right)  $
permutations of the form $s_{k}$ (with $k\in\left\{  1,2,\ldots,n-1\right\}
$). In other words, there exists an $\ell\left(  \tau\right)  $-tuple $\left(
k_{1},k_{2},\ldots,k_{\ell\left(  \tau\right)  }\right)  \in\left\{
1,2,\ldots,n-1\right\}  ^{\ell\left(  \tau\right)  }$ such that $\tau
=s_{k_{1}}\circ s_{k_{2}}\circ\cdots\circ s_{k_{\ell\left(  \tau\right)  }}$.
Consider this $\left(  k_{1},k_{2},\ldots,k_{\ell\left(  \tau\right)
}\right)  $. Then,%
\[
\ell\left(  \sigma\circ\underbrace{\tau}_{=s_{k_{1}}\circ s_{k_{2}}\circ
\cdots\circ s_{k_{\ell\left(  \tau\right)  }}}\right)  =\ell\left(
\sigma\circ\left(  s_{k_{1}}\circ s_{k_{2}}\circ\cdots\circ s_{k_{\ell\left(
\tau\right)  }}\right)  \right)  \leq\ell\left(  \sigma\right)  +\ell\left(
\tau\right)
\]
(by (\ref{sol.ps2.2.5.c.snake}), applied to $p=\ell\left(  \tau\right)  $).
This solves Exercise \ref{exe.ps2.2.5} \textbf{(c)}.

\textbf{(g)} Exercise \ref{exe.ps2.2.5} \textbf{(e)} shows that $\sigma$ can
be written as a composition of $\ell\left(  \sigma\right)  $ permutations of
the form $s_{k}$ (with $k\in\left\{  1,2,\ldots,n-1\right\}  $). In other
words, $\ell\left(  \sigma\right)  $ is an $N\in\mathbb{N}$ such that $\sigma$
can be written as a composition of $N$ permutations of the form $s_{k}$ (with
$k\in\left\{  1,2,\ldots,n-1\right\}  $). In order to solve Exercise
\ref{exe.ps2.2.5} \textbf{(g)}, it only remains to show that $\ell\left(
\sigma\right)  $ is the \textbf{smallest} such $N$. In other words, it remains
to show that if $N\in\mathbb{N}$ is such that $\sigma$ can be written as a
composition of $N$ permutations of the form $s_{k}$ (with $k\in\left\{
1,2,\ldots,n-1\right\}  $), then $N\geq\ell\left(  \sigma\right)  $.

So let $N\in\mathbb{N}$ be such that $\sigma$ can be written as a composition
of $N$ permutations of the form $s_{k}$ (with $k\in\left\{  1,2,\ldots
,n-1\right\}  $). In other words, there exists an $N$-tuple $\left(
k_{1},k_{2},\ldots,k_{N}\right)  \in\left\{  1,2,\ldots,n-1\right\}  ^{N}$
such that $\sigma=s_{k_{1}}\circ s_{k_{2}}\circ\cdots\circ s_{k_{N}}$.
Applying (\ref{sol.ps2.2.5.c.snake}) to $\operatorname*{id}$ and $N$ instead
of $\sigma$ and $p$, we obtain%
\[
\ell\left(  \operatorname*{id}\circ\left(  s_{k_{1}}\circ s_{k_{2}}\circ
\cdots\circ s_{k_{N}}\right)  \right)  \leq\underbrace{\ell\left(
\operatorname*{id}\right)  }_{=0}+N=N.
\]
Since $\operatorname*{id}\circ\left(  s_{k_{1}}\circ s_{k_{2}}\circ\cdots\circ
s_{k_{N}}\right)  =s_{k_{1}}\circ s_{k_{2}}\circ\cdots\circ s_{k_{N}}=\sigma$,
this rewrites as $\ell\left(  \sigma\right)  \leq N$. In other words,
$N\geq\ell\left(  \sigma\right)  $. This completes our solution to Exercise
\ref{exe.ps2.2.5} \textbf{(g)}.
\end{proof}

\begin{remark}
The above solution to Exercise \ref{exe.ps2.2.5} owes most of its length to my
attempts at being precise. As Pascal said, \textquotedblleft I have made this
longer than usual because I have not had time to make it
shorter\textquotedblright. The proofs are not, per se, difficult, but this is
combinatorics, and proofs in combinatorics often have to walk a tightrope
between being unreadably long and unreadably terse.

Most parts of Exercise \ref{exe.ps2.2.5} can be proven in more than just one
way. Let me briefly mention an alternative proof for parts \textbf{(b)} and
\textbf{(c)}. Namely, let us use the notation $\operatorname*{Inv}\left(
\sigma\right)  $ for the set of all inversions of a permutation $\sigma$. Let
$\sigma\in S_{n}$ and $\tau\in S_{n}$. Then, it is not hard to convince
oneself that%
\[
\operatorname*{Inv}\left(  \sigma\circ\tau\right)  =A\cup B,
\]
where $A$ and $B$ are the sets defined by%
\begin{align*}
A  &  =\left\{  \left(  i,j\right)  \ \mid\ \left(  i,j\right)  \in
\operatorname*{Inv}\left(  \tau\right)  \text{ and }\left(  \tau\left(
j\right)  ,\tau\left(  i\right)  \right)  \notin\operatorname*{Inv}\left(
\sigma\right)  \right\}  ;\\
B  &  =\left\{  \left(  i,j\right)  \ \mid\ \left(  i,j\right)  \notin%
\operatorname*{Inv}\left(  \tau\right)  \text{ and }\left(  \tau\left(
i\right)  ,\tau\left(  j\right)  \right)  \in\operatorname*{Inv}\left(
\sigma\right)  \right\}
\end{align*}
(where $\left(  i,j\right)  $ are subject to the condition $1\leq i<j\leq n$
both times). (Essentially, this is because an $\left(  i,j\right)
\in\operatorname*{Inv}\left(  \sigma\circ\tau\right)  $ either satisfies
$\tau\left(  i\right)  >\tau\left(  j\right)  $ or satisfies $\tau\left(
i\right)  <\tau\left(  j\right)  $. In the first case, this $\left(
i,j\right)  $ belongs to $A$; in the second case, it belongs to $B$.) It is
furthermore clear that $\left\vert A\right\vert \leq\left\vert
\operatorname*{Inv}\left(  \tau\right)  \right\vert =\ell\left(  \tau\right)
$ (by (\ref{sol.ps2.2.5.f.1})) and $\left\vert B\right\vert \leq\left\vert
\operatorname*{Inv}\left(  \sigma\right)  \right\vert =\ell\left(
\sigma\right)  $ (by (\ref{sol.ps2.2.5.f.1})). Hence, (\ref{sol.ps2.2.5.f.1})
yields $\ell\left(  \sigma\circ\tau\right)  =\left\vert
\underbrace{\operatorname*{Inv}\left(  \sigma\circ\tau\right)  }_{=A\cup
B}\right\vert =\left\vert A\cup B\right\vert \leq\underbrace{\left\vert
A\right\vert }_{\leq\ell\left(  \tau\right)  }+\underbrace{\left\vert
B\right\vert }_{\leq\ell\left(  \sigma\right)  }\leq\ell\left(  \tau\right)
+\ell\left(  \sigma\right)  $. This solves Exercise \ref{exe.ps2.2.5}
\textbf{(c)}. To solve part \textbf{(b)}, we need to take a few more steps.
First, it is clear that $A\cap B=\varnothing$, so that $\left\vert A\cup
B\right\vert =\left\vert A\right\vert +\left\vert B\right\vert $. Second, let
us set%
\[
C=\left\{  \left(  i,j\right)  \in\left\{  1,2,\ldots,n\right\}  ^{2}%
\ \mid\ i<j\text{, }\tau\left(  i\right)  >\tau\left(  j\right)  \text{ and
}\sigma\left(  \tau\left(  i\right)  \right)  <\sigma\left(  \tau\left(
j\right)  \right)  \right\}  .
\]
Then, it is easy to see that $C\subseteq\operatorname*{Inv}\left(
\tau\right)  $ and $A=\operatorname*{Inv}\left(  \tau\right)  \setminus C$, so
that $\left\vert A\right\vert =\left\vert \operatorname*{Inv}\left(
\tau\right)  \right\vert -\left\vert C\right\vert $. Moreover, if $\mathbf{t}$
denotes the permutation of $\left\{  1,2,\ldots,n\right\}  ^{2}$ which sends
every $\left(  i,j\right)  \in\left\{  1,2,\ldots,n\right\}  ^{2}$ to $\left(
\tau\left(  i\right)  ,\tau\left(  j\right)  \right)  $, and if $\mathbf{f}$
denotes the permutation of $\left\{  1,2,\ldots,n\right\}  ^{2}$ which sends
every $\left(  i,j\right)  $ to $\left(  j,i\right)  $, then $\mathbf{f}%
\left(  C\right)  \subseteq\mathbf{t}^{-1}\left(  \operatorname*{Inv}\left(
\sigma\right)  \right)  $ and $B=\mathbf{t}^{-1}\left(  \operatorname*{Inv}%
\left(  \sigma\right)  \right)  \setminus\mathbf{f}\left(  C\right)  $, so
that $\left\vert B\right\vert =\underbrace{\left\vert \mathbf{t}^{-1}\left(
\operatorname*{Inv}\left(  \sigma\right)  \right)  \right\vert }_{=\left\vert
\operatorname*{Inv}\left(  \sigma\right)  \right\vert }-\underbrace{\left\vert
\mathbf{f}\left(  C\right)  \right\vert }_{=\left\vert C\right\vert
}=\left\vert \operatorname*{Inv}\left(  \sigma\right)  \right\vert -\left\vert
C\right\vert $. Thus,%
\begin{align*}
\ell\left(  \sigma\circ\tau\right)   &  =\left\vert
\underbrace{\operatorname*{Inv}\left(  \sigma\circ\tau\right)  }_{=A\cup
B}\right\vert =\left\vert A\cup B\right\vert =\underbrace{\left\vert
A\right\vert }_{=\left\vert \operatorname*{Inv}\left(  \tau\right)
\right\vert -\left\vert C\right\vert }+\underbrace{\left\vert B\right\vert
}_{=\left\vert \operatorname*{Inv}\left(  \sigma\right)  \right\vert
-\left\vert C\right\vert }\\
&  =\left(  \underbrace{\left\vert \operatorname*{Inv}\left(  \tau\right)
\right\vert }_{=\ell\left(  \tau\right)  }-\left\vert C\right\vert \right)
+\left(  \underbrace{\left\vert \operatorname*{Inv}\left(  \sigma\right)
\right\vert }_{=\ell\left(  \sigma\right)  }-\left\vert C\right\vert \right)
=\left(  \ell\left(  \tau\right)  -\left\vert C\right\vert \right)  +\left(
\ell\left(  \sigma\right)  -\left\vert C\right\vert \right) \\
&  =\ell\left(  \sigma\right)  +\ell\left(  \tau\right)  -2\left\vert
C\right\vert \equiv\ell\left(  \sigma\right)  +\ell\left(  \tau\right)
\operatorname{mod}2.
\end{align*}
This solves Exercise \ref{exe.ps2.2.5} \textbf{(b)}.
\end{remark}

\subsection{Solution to Exercise \ref{exe.ps2.2.6}}

\begin{proof}
[Solution to Exercise \ref{exe.ps2.2.6}.]We need to show that if
$p\in\mathbb{N}$ and \newline$\left(  k_{1},k_{2},\ldots,k_{p}\right)
\in\left\{  1,2,\ldots,n-1\right\}  ^{p}$ are such that $\sigma=s_{k_{1}}\circ
s_{k_{2}}\circ\cdots\circ s_{k_{p}}$, then $p\equiv\ell\left(  \sigma\right)
\operatorname{mod}2$.

Let $p\in\mathbb{N}$ and $\left(  k_{1},k_{2},\ldots,k_{p}\right)  \in\left\{
1,2,\ldots,n-1\right\}  ^{p}$ be such that $\sigma=s_{k_{1}}\circ s_{k_{2}%
}\circ\cdots\circ s_{k_{p}}$. We must prove that $p\equiv\ell\left(
\sigma\right)  \operatorname{mod}2$.

Applying (\ref{sol.ps2.2.5.b.snake}) to $\operatorname*{id}$ instead of
$\sigma$, we see that%
\[
\ell\left(  \operatorname*{id}\circ\left(  s_{k_{1}}\circ s_{k_{2}}\circ
\cdots\circ s_{k_{p}}\right)  \right)  \equiv\underbrace{\ell\left(
\operatorname*{id}\right)  }_{=0}+p=p\operatorname{mod}2.
\]
Since $\operatorname*{id}\circ\left(  s_{k_{1}}\circ s_{k_{2}}\circ\cdots\circ
s_{k_{p}}\right)  =s_{k_{1}}\circ s_{k_{2}}\circ\cdots\circ s_{k_{p}}=\sigma$,
this rewrites as $\ell\left(  \sigma\right)  \equiv p\operatorname{mod}2$. In
other words, $p\equiv\ell\left(  \sigma\right)  \operatorname{mod}2$. This
completes the solution of Exercise \ref{exe.ps2.2.6}.
\end{proof}

\subsection{Solution to Exercise \ref{exe.ps2.2.7}}

\begin{proof}
[Solution to Exercise \ref{exe.ps2.2.7}.]We have $n\geq2$. Thus, the
permutation $s_{1}\in S_{n}$ is well-defined. (This is the permutation which
swaps $1$ with $2$ while leaving all other elements of $\left\{
1,2,\ldots,n\right\}  $ unchanged.)

Let $A_{n}$ denote the set of all even permutations in $S_{n}$. Let $C_{n}$
denote the set of all odd permutations in $S_{n}$. The sign of a permutation
in $S_{n}$ is either $1$ or $-1$ (because it is defined to be an integer power
of $-1$), but not both. Hence, every permutation in $S_{n}$ is either even or
odd, but not both. In other words, we have $S_{n}=A_{n}\cup C_{n}$ and
$A_{n}\cap C_{n}=\varnothing$. Therefore, $\left\vert S_{n}\right\vert
=\left\vert A_{n}\right\vert +\left\vert C_{n}\right\vert $.

Now, for every $\sigma\in A_{n}$, we have $\sigma\circ s_{k}\in C_{n}%
$\ \ \ \ \footnote{\textit{Proof.} Let $\sigma\in A_{n}$. Thus, $\sigma$ is an
even permutation in $S_{n}$ (since $A_{n}$ is the set of all even permutations
in $S_{n}$). Since $\sigma$ is even, we have $\left(  -1\right)  ^{\sigma}=1$,
so that $1=\left(  -1\right)  ^{\sigma}=\left(  -1\right)  ^{\ell\left(
\sigma\right)  }$. Therefore, $\ell\left(  \sigma\right)  \equiv
0\operatorname{mod}2$. Now, (\ref{sol.ps2.2.5.b.stepper}) yields $\ell\left(
\sigma\circ s_{k}\right)  \equiv\underbrace{\ell\left(  \sigma\right)
}_{\equiv0\operatorname{mod}2}+1\equiv1\operatorname{mod}2$, so that $\left(
-1\right)  ^{\ell\left(  \sigma\circ s_{k}\right)  }=-1$. But now, $\left(
-1\right)  ^{\sigma\circ s_{k}}=\left(  -1\right)  ^{\ell\left(  \sigma\circ
s_{k}\right)  }=-1$, so that the permutation $\sigma\circ s_{k}$ is odd. In
other words, $\sigma\circ s_{k}\in C_{n}$ (since $C_{n}$ is the set of all odd
permutations in $S_{n}$), qed.}. Hence, we can define a map $\Phi
:A_{n}\rightarrow C_{n}$ by
\[
\Phi\left(  \sigma\right)  =\sigma\circ s_{k}\ \ \ \ \ \ \ \ \ \ \text{for
every }\sigma\in A_{n}.
\]
Similarly, we can define a map $\Psi:C_{n}\rightarrow A_{n}$ by
\[
\Psi\left(  \sigma\right)  =\sigma\circ s_{k}\ \ \ \ \ \ \ \ \ \ \text{for
every }\sigma\in C_{n}.
\]
These two maps $\Phi$ and $\Psi$ are mutually inverse (since every $\sigma\in
S_{n}$ satisfies $\sigma\circ\underbrace{s_{k}\circ s_{k}}_{=s_{k}%
^{2}=\operatorname*{id}}=\sigma$). Therefore, the map $\Phi$ is a bijection.
Thus, there exists a bijection form $A_{n}$ to $C_{n}$ (namely, $\Phi$), so
that we obtain $\left\vert A_{n}\right\vert =\left\vert C_{n}\right\vert $.
Hence, $\left\vert S_{n}\right\vert =\underbrace{\left\vert A_{n}\right\vert
}_{=\left\vert C_{n}\right\vert }+\left\vert C_{n}\right\vert =\left\vert
C_{n}\right\vert +\left\vert C_{n}\right\vert =2\left\vert C_{n}\right\vert $
and therefore $\left\vert C_{n}\right\vert =\dfrac{1}{2}\underbrace{\left\vert
S_{n}\right\vert }_{=n!}=\dfrac{1}{2}n!=n!/2$. In other words, the number of
odd permutations in $S_{n}$ is $n!/2$. Similarly, the number of even
permutations in $S_{n}$ is $n!/2$. Exercise \ref{exe.ps2.2.7} is solved.
\end{proof}

\subsection{Solution to Exercise \ref{exe.ps2.2.4'}}

\begin{proof}
[Solution to Exercise \ref{exe.ps2.2.4'}.]The solution of Exercise
\ref{exe.ps2.2.4'} \textbf{(a)} is completely analogous to the solution of
Exercise \ref{exe.ps2.2.4} \textbf{(a)}; it can be obtained from the latter by
replacing $S_{n}$ by $S_{\infty}$, replacing $\left\{  1,2,\ldots,n\right\}  $
by $\left\{  1,2,3,\ldots\right\}  $, and replacing $\left\{  1,2,\ldots
,n-2\right\}  $ by $\left\{  1,2,3,\ldots\right\}  $.

As for Exercise \ref{exe.ps2.2.4'} \textbf{(b)}, we again omit the solution,
because it follows from an exercise that will be solved below (Exercise
\ref{exe.ps2.2.5'} \textbf{(e)}).
\end{proof}

\subsection{Solution to Exercise \ref{exe.ps2.2.5'}}

\begin{proof}
[Solution to Exercise \ref{exe.ps2.2.5'}.]A solution of Exercise
\ref{exe.ps2.2.5'} can be obtained by copying our above solution of Exercise
\ref{exe.ps2.2.5} almost verbatim, occasionally doing some replacements (e.g.,
we have to replace $S_{n}$ by $S_{\left(  \infty\right)  }$, to replace $1\leq
i<j\leq n$ by $1\leq i<j$, to replace $\left\{  1,2,\ldots,n-1\right\}  $ by
$\left\{  1,2,3,\ldots\right\}  $, and to replace $\left\{  1,2,\ldots
,n\right\}  $ by $\left\{  1,2,3,\ldots\right\}  $). We leave the
straightforward changes to the reader.
\end{proof}

\subsection{Solution to Exercise \ref{exe.ps2.2.6'}}

\begin{proof}
[Solution to Exercise \ref{exe.ps2.2.6'}.]A solution of Exercise
\ref{exe.ps2.2.6'} can be obtained by copying our above solution of Exercise
\ref{exe.ps2.2.6} almost verbatim, occasionally doing some replacements (e.g.,
we have to replace $S_{n}$ by $S_{\left(  \infty\right)  }$, to replace
$\left\{  1,2,\ldots,n-1\right\}  $ by $\left\{  1,2,3,\ldots\right\}  $, and
to replace $\left\{  1,2,\ldots,n\right\}  $ by $\left\{  1,2,3,\ldots
\right\}  $). We leave the straightforward changes to the reader.
\end{proof}

\subsection{Solution to Exercise \ref{exe.ps4.0}}

\begin{vershort}
\begin{proof}
[Solution to Exercise \ref{exe.ps4.0}.]In the following, \textquotedblleft
path\textquotedblright\ will always mean \textquotedblleft path on the
(undirected) $n$-th right Bruhat graph\textquotedblright. Hence, we need to
prove that $\ell\left(  \sigma^{-1}\circ\tau\right)  $ is the smallest length
of a path between $\sigma$ and $\tau$.

We write any path as the tuple consisting of its vertices (from its beginning
to its end).\footnote{This is legitimate, because the (undirected) $n$-th
right Bruhat graph does not have multiple edges.}

Let $L=\ell\left(  \sigma^{-1}\circ\tau\right)  $. We shall first show that
there exists a path of length $L$ between $\sigma$ and $\tau$.

Indeed, Exercise \ref{exe.ps2.2.5} \textbf{(g)} (applied to $\sigma^{-1}%
\circ\tau$ instead of $\sigma$) yields that $\ell\left(  \sigma^{-1}\circ
\tau\right)  $ is the smallest $N\in\mathbb{N}$ such that $\sigma^{-1}%
\circ\tau$ can be written as a composition of $N$ permutations of the form
$s_{k}$ (with $k\in\left\{  1,2,\ldots,n-1\right\}  $). Since $L=\ell\left(
\sigma^{-1}\circ\tau\right)  $, we can rewrite this as follows: $L$ is the
smallest $N\in\mathbb{N}$ such that $\sigma^{-1}\circ\tau$ can be written as a
composition of $N$ permutations of the form $s_{k}$ (with $k\in\left\{
1,2,\ldots,n-1\right\}  $). Hence, $\sigma^{-1}\circ\tau$ can be written as a
composition of $L$ permutations of the form $s_{k}$ (with $k\in\left\{
1,2,\ldots,n-1\right\}  $). In other words, there exists an $L$-tuple $\left(
j_{1},j_{2},\ldots,j_{L}\right)  \in\left\{  1,2,\ldots,n-1\right\}  ^{L}$
such that $\sigma^{-1}\circ\tau=s_{j_{1}}\circ s_{j_{2}}\circ\cdots\circ
s_{j_{L}}$. Consider this $\left(  j_{1},j_{2},\ldots,j_{L}\right)  $.

We have $\underbrace{\sigma\circ\sigma^{-1}}_{=\operatorname*{id}}\circ
\tau=\tau$, so that $\tau=\sigma\circ\underbrace{\sigma^{-1}\circ\tau
}_{=s_{j_{1}}\circ s_{j_{2}}\circ\cdots\circ s_{j_{L}}}=\sigma\circ\left(
s_{j_{1}}\circ s_{j_{2}}\circ\cdots\circ s_{j_{L}}\right)  $. Now, for every
$p\in\left\{  0,1,\ldots,L\right\}  $, define a permutation $\gamma_{p}\in
S_{n}$ by%
\[
\gamma_{p}=\sigma\circ\left(  s_{j_{1}}\circ s_{j_{2}}\circ\cdots\circ
s_{j_{p}}\right)  .
\]
Thus, $\gamma_{0}=\sigma\circ\underbrace{\left(  s_{j_{1}}\circ s_{j_{2}}%
\circ\cdots\circ s_{j_{0}}\right)  }_{=\left(  \text{a composition of }0\text{
permutations}\right)  =\operatorname*{id}}=\sigma$ and $\gamma_{L}=\sigma
\circ\left(  s_{j_{1}}\circ s_{j_{2}}\circ\cdots\circ s_{j_{L}}\right)  =\tau$.

Now, for every $i\in\left\{  1,2,\ldots,L\right\}  $, the vertices $\gamma
_{i}$ and $\gamma_{i-1}$ of the (undirected) $n$-th right Bruhat graph are
adjacent\footnote{\textit{Proof.} Let $i\in\left\{  1,2,\ldots,L\right\}  $.
Then, the definition of $\gamma_{i-1}$ yields $\gamma_{i-1}=\sigma\circ\left(
s_{j_{1}}\circ s_{j_{2}}\circ\cdots\circ s_{j_{i-1}}\right)  $, whereas the
definition of $\gamma_{i}$ yields%
\[
\gamma_{i}=\sigma\circ\underbrace{\left(  s_{j_{1}}\circ s_{j_{2}}\circ
\cdots\circ s_{j_{i}}\right)  }_{=\left(  s_{j_{1}}\circ s_{j_{2}}\circ
\cdots\circ s_{j_{i-1}}\right)  \circ s_{j_{i}}}=\underbrace{\sigma
\circ\left(  s_{j_{1}}\circ s_{j_{2}}\circ\cdots\circ s_{j_{i-1}}\right)
}_{=\gamma_{i-1}}\circ s_{j_{i}}=\gamma_{i-1}\circ s_{j_{i}}.
\]
Therefore, there exists a $k\in\left\{  1,2,\ldots,n-1\right\}  $ such that
$\gamma_{i}=\gamma_{i-1}\circ s_{k}$ (namely, $k=j_{i}$). In other words, the
vertices $\gamma_{i}$ and $\gamma_{i-1}$ of the (undirected) $n$-th right
Bruhat graph are adjacent (by the definition of the edges of this graph).
Qed.}. Hence, the vertices $\gamma_{0},\gamma_{1},\ldots,\gamma_{L}$ form a
path. This path connects $\sigma$ to $\tau$ (since it begins at $\gamma
_{0}=\sigma$ and ends at $\gamma_{L}=\tau$), and has length $L$. Thus, there
exists a path between $\sigma$ and $\tau$ which has length $L$ (namely, the
path formed by the vertices $\gamma_{0},\gamma_{1},\ldots,\gamma_{L}$).

We shall now show that $L$ is the smallest length of a path between $\sigma$
and $\tau$. Indeed, let $\mathbf{d}$ be any path between $\sigma$ and $\tau$.
We shall show that the length of $\mathbf{d}$ is $\geq L$.

The path $\mathbf{d}$ is a path between $\sigma$ and $\tau$. Hence, we can
write the path $\mathbf{d}$ in the form $\mathbf{d}=\left(  \delta_{0}%
,\delta_{1},\ldots,\delta_{M}\right)  $ for some $\delta_{0},\delta_{1}%
,\ldots,\delta_{M}\in S_{n}$ with $\delta_{0}=\sigma$ and $\delta_{M}=\tau$.
Consider these $\delta_{0},\delta_{1},\ldots,\delta_{M}\in S_{n}$.

For every $i\in\left\{  1,2,\ldots,M\right\}  $, there exists a $k\in\left\{
1,2,\ldots,n-1\right\}  $ such that $\delta_{i}=\delta_{i-1}\circ s_{k}%
$\ \ \ \ \footnote{\textit{Proof.} Let $i\in\left\{  1,2,\ldots,M\right\}  $.
Then, the vertices $\delta_{i}$ and $\delta_{i-1}$ of the (undirected) $n$-th
right Bruhat graph are adjacent (because they are two consecutive vertices on
the path $\left(  \delta_{0},\delta_{1},\ldots,\delta_{M}\right)  =\mathbf{d}%
$). In other words, there exists a $k\in\left\{  1,2,\ldots,n-1\right\}  $
such that $\delta_{i}=\delta_{i-1}\circ s_{k}$ (by the definition of the edges
of this graph). Qed.}. We denote this $k$ by $k_{i}$. Thus, for every
$i\in\left\{  1,2,\ldots,M\right\}  $, we have defined a $k_{i}\in\left\{
1,2,\ldots,n-1\right\}  $ such that $\delta_{i}=\delta_{i-1}\circ s_{k_{i}}$.

Now, every $j\in\left\{  0,1,\ldots,M\right\}  $ satisfies%
\begin{equation}
\delta_{j}=\sigma\circ\left(  s_{k_{1}}\circ s_{k_{2}}\circ\cdots\circ
s_{k_{j}}\right)  \label{sol.ps4.short.0.5}%
\end{equation}
\footnote{\textit{Proof of (\ref{sol.ps4.short.0.5}):} We shall prove
(\ref{sol.ps4.short.0.5}) by induction over $j$:
\par
\textit{Induction base:} We have $\delta_{0}=\sigma$. Compared with
$\sigma\circ\underbrace{\left(  s_{k_{1}}\circ s_{k_{2}}\circ\cdots\circ
s_{k_{0}}\right)  }_{=\left(  \text{a composition of }0\text{ permutations}%
\right)  =\operatorname*{id}}=\sigma$, this yields $\delta_{0}=\sigma
\circ\left(  s_{k_{1}}\circ s_{k_{2}}\circ\cdots\circ s_{k_{0}}\right)  $. In
other words, (\ref{sol.ps4.short.0.5}) holds for $j=0$. This completes the
induction base.
\par
\textit{Induction step:} Let $J\in\left\{  0,1,\ldots,M\right\}  $ be
positive. Assume that (\ref{sol.ps4.short.0.5}) holds for $j=J-1$. We need to
show that (\ref{sol.ps4.short.0.5}) holds for $j=J$.
\par
We have $\delta_{J}=\sigma\circ\left(  s_{k_{1}}\circ s_{k_{2}}\circ
\cdots\circ s_{k_{J-1}}\right)  $ (since (\ref{sol.ps4.short.0.5}) holds for
$j=J-1$). Now, recall that $\delta_{i}=\delta_{i-1}\circ s_{k_{i}}$ for every
$i\in\left\{  1,2,\ldots,M\right\}  $. Applying this to $i=J$, we obtain%
\begin{align*}
\delta_{J}  &  =\underbrace{\delta_{J-1}}_{=\sigma\circ\left(  s_{k_{1}}\circ
s_{k_{2}}\circ\cdots\circ s_{k_{J-1}}\right)  }\circ s_{k_{J}}%
\ \ \ \ \ \ \ \ \ \ \left(  \text{since }J\in\left\{  1,2,\ldots,M\right\}
\text{ (since }J\in\left\{  0,1,\ldots,M\right\}  \text{ is positive)}\right)
\\
&  =\sigma\circ\underbrace{\left(  s_{k_{1}}\circ s_{k_{2}}\circ\cdots\circ
s_{k_{J-1}}\right)  \circ s_{k_{J}}}_{=s_{k_{1}}\circ s_{k_{2}}\circ
\cdots\circ s_{k_{J}}}=\sigma\circ\left(  s_{k_{1}}\circ s_{k_{2}}\circ
\cdots\circ s_{k_{J}}\right)  .
\end{align*}
In other words, (\ref{sol.ps4.short.0.5}) holds for $j=J$. This completes the
induction step. Thus, (\ref{sol.ps4.short.0.5}) is proven by induction.}.
Applying this to $j=M$, we obtain%
\[
\delta_{M}=\sigma\circ\left(  s_{k_{1}}\circ s_{k_{2}}\circ\cdots\circ
s_{k_{M}}\right)  .
\]
Compared with $\delta_{M}=\tau$, this yields%
\[
\tau=\sigma\circ\left(  s_{k_{1}}\circ s_{k_{2}}\circ\cdots\circ s_{k_{M}%
}\right)  ,
\]
so that
\[
\sigma^{-1}\circ\underbrace{\tau}_{=\sigma\circ\left(  s_{k_{1}}\circ
s_{k_{2}}\circ\cdots\circ s_{k_{M}}\right)  }=\underbrace{\sigma^{-1}%
\circ\sigma}_{=\operatorname*{id}}\circ\left(  s_{k_{1}}\circ s_{k_{2}}%
\circ\cdots\circ s_{k_{M}}\right)  =s_{k_{1}}\circ s_{k_{2}}\circ\cdots\circ
s_{k_{M}}.
\]
Therefore, the permutation $\sigma^{-1}\circ\tau$ can be written as a
composition of $M$ permutations of the form $s_{k}$ (with $k\in\left\{
1,2,\ldots,n-1\right\}  $) (namely, of the $M$ permutations $s_{k_{1}}$,
$s_{k_{2}}$, $\ldots$, $s_{k_{M}}$).

Now, we recall that $L$ is the smallest $N\in\mathbb{N}$ such that
$\sigma^{-1}\circ\tau$ can be written as a composition of $N$ permutations of
the form $s_{k}$ (with $k\in\left\{  1,2,\ldots,n-1\right\}  $). Hence, if
$N\in\mathbb{N}$ is such that $\sigma^{-1}\circ\tau$ can be written as a
composition of $N$ permutations of the form $s_{k}$ (with $k\in\left\{
1,2,\ldots,n-1\right\}  $), then $N\geq L$. We can apply this to $N=M$
(because $\sigma^{-1}\circ\tau$ can be written as a composition of $M$
permutations of the form $s_{k}$ (with $k\in\left\{  1,2,\ldots,n-1\right\}
$)), and thus obtain $M\geq L$.

But the length of the path $\mathbf{d}$ is $M$ (since $\mathbf{d}=\left(
\delta_{0},\delta_{1},\ldots,\delta_{M}\right)  $). Hence, the length of the
path $\mathbf{d}$ is $\geq L$ (since $M\geq L$).

Let us now forget that we fixed $\mathbf{d}$. We thus have shown that if
$\mathbf{d}$ is any path between $\sigma$ and $\tau$, then the length of the
path $\mathbf{d}$ is $\geq L$. In other words, every path between $\sigma$ and
$\tau$ has length $\geq L$.

Altogether, we have proven the following two statements:

\begin{itemize}
\item There exists a path of length $L$ between $\sigma$ and $\tau$.

\item Every path between $\sigma$ and $\tau$ has length $\geq L$.
\end{itemize}

Therefore, $L$ is the smallest length of a path between $\sigma$ and $\tau$.
In other words, $\ell\left(  \sigma^{-1}\circ\tau\right)  $ is the smallest
length of a path between $\sigma$ and $\tau$ (since $L=\ell\left(  \sigma
^{-1}\circ\tau\right)  $). Exercise \ref{exe.ps4.0} is solved.
\end{proof}
\end{vershort}

\begin{verlong}
\begin{proof}
[Solution to Exercise \ref{exe.ps4.0}.]Let us first forget about $\sigma$ and
$\tau$, and make a few general remarks.

Let $\sigma\in S_{n}$. Then, Exercise \ref{exe.ps2.2.5} \textbf{(g)} shows
that
\begin{equation}
\left(
\begin{array}
[c]{c}%
\ell\left(  \sigma\right)  \text{ is the smallest }N\in\mathbb{N}\text{ such
that }\sigma\text{ can be written as a composition}\\
\text{of }N\text{ permutations of the form }s_{k}\text{ (with }k\in\left\{
1,2,\ldots,n-1\right\}  \text{)}%
\end{array}
\right)  . \label{sol.ps4.0.lem}%
\end{equation}

Let us now forget that we fixed $\sigma$. We thus have proven
(\ref{sol.ps4.0.lem}) for every $\sigma\in S_{n}$.

Now, let $\sigma\in S_{n}$ and $\tau\in S_{n}$. In the following,
\textquotedblleft path\textquotedblright\ will always mean \textquotedblleft
path on the (undirected) $n$-th right Bruhat graph\textquotedblright. Hence,
we need to prove that $\ell\left(  \sigma^{-1}\circ\tau\right)  $ is the
smallest length of a path between $\sigma$ and $\tau$.

We write any path as the tuple consisting of its vertices (from its beginning
to its end).\footnote{This is legitimate, because the (undirected) $n$-th
right Bruhat graph does not have multiple edges.}

Let $L=\ell\left(  \sigma^{-1}\circ\tau\right)  $. We shall first show that
there exists a path of length $L$ between $\sigma$ and $\tau$.

Indeed, (\ref{sol.ps4.0.lem}) (applied to $\sigma^{-1}\circ\tau$ instead of
$\sigma$) yields that $\ell\left(  \sigma^{-1}\circ\tau\right)  $ is the
smallest $N\in\mathbb{N}$ such that $\sigma^{-1}\circ\tau$ can be written as a
composition of $N$ permutations of the form $s_{k}$ (with $k\in\left\{
1,2,\ldots,n-1\right\}  $). Hence, $\ell\left(  \sigma^{-1}\circ\tau\right)  $
is an $N\in\mathbb{N}$ such that $\sigma^{-1}\circ\tau$ can be written as a
composition of $N$ permutations of the form $s_{k}$ (with $k\in\left\{
1,2,\ldots,n-1\right\}  $). In other words, $\sigma^{-1}\circ\tau$ can be
written as a composition of $\ell\left(  \sigma^{-1}\circ\tau\right)  $
permutations of the form $s_{k}$ (with $k\in\left\{  1,2,\ldots,n-1\right\}
$). In other words, $\sigma^{-1}\circ\tau$ can be written as a composition of
$L$ permutations of the form $s_{k}$ (with $k\in\left\{  1,2,\ldots
,n-1\right\}  $) (since $L=\ell\left(  \sigma^{-1}\circ\tau\right)  $). In
other words, there exists an $L$-tuple $\left(  j_{1},j_{2},\ldots
,j_{L}\right)  \in\left\{  1,2,\ldots,n-1\right\}  ^{L}$ such that
$\sigma^{-1}\circ\tau=s_{j_{1}}\circ s_{j_{2}}\circ\cdots\circ s_{j_{L}}$.
Consider this $\left(  j_{1},j_{2},\ldots,j_{L}\right)  $.

We have $\underbrace{\sigma\circ\sigma^{-1}}_{=\operatorname*{id}}\circ
\tau=\tau$, so that $\tau=\sigma\circ\underbrace{\sigma^{-1}\circ\tau
}_{=s_{j_{1}}\circ s_{j_{2}}\circ\cdots\circ s_{j_{L}}}=\sigma\circ\left(
s_{j_{1}}\circ s_{j_{2}}\circ\cdots\circ s_{j_{L}}\right)  $. Now, for every
$p\in\left\{  0,1,\ldots,L\right\}  $, define a permutation $\gamma_{p}\in
S_{n}$ by%
\[
\gamma_{p}=\sigma\circ\left(  s_{j_{1}}\circ s_{j_{2}}\circ\cdots\circ
s_{j_{p}}\right)  .
\]
Thus, $\gamma_{0}=\sigma\circ\underbrace{\left(  s_{j_{1}}\circ s_{j_{2}}%
\circ\cdots\circ s_{j_{0}}\right)  }_{=\left(  \text{a composition of }0\text{
permutations}\right)  =\operatorname*{id}}=\sigma$ and $\gamma_{L}=\sigma
\circ\left(  s_{j_{1}}\circ s_{j_{2}}\circ\cdots\circ s_{j_{L}}\right)  =\tau$.

Now, for every $i\in\left\{  1,2,\ldots,L\right\}  $, the vertices $\gamma
_{i}$ and $\gamma_{i-1}$ of the (undirected) $n$-th right Bruhat graph are
adjacent\footnote{\textit{Proof.} Let $i\in\left\{  1,2,\ldots,L\right\}  $.
Then, the definition of $\gamma_{i-1}$ yields $\gamma_{i-1}=\sigma\circ\left(
s_{j_{1}}\circ s_{j_{2}}\circ\cdots\circ s_{j_{i-1}}\right)  $, whereas the
definition of $\gamma_{i}$ yields%
\begin{align*}
\gamma_{i}  &  =\sigma\circ\underbrace{\left(  s_{j_{1}}\circ s_{j_{2}}%
\circ\cdots\circ s_{j_{i}}\right)  }_{=\left(  s_{j_{1}}\circ s_{j_{2}}%
\circ\cdots\circ s_{j_{i-1}}\right)  \circ s_{j_{i}}}=\underbrace{\sigma
\circ\left(  s_{j_{1}}\circ s_{j_{2}}\circ\cdots\circ s_{j_{i-1}}\right)
}_{=\gamma_{i-1}}\circ s_{j_{i}}\\
&  =\gamma_{i-1}\circ s_{j_{i}}.
\end{align*}
Therefore, there exists a $k\in\left\{  1,2,\ldots,n-1\right\}  $ such that
$\gamma_{i}=\gamma_{i-1}\circ s_{k}$ (namely, $k=j_{i}$). In other words, the
vertices $\gamma_{i}$ and $\gamma_{i-1}$ of the (undirected) $n$-th right
Bruhat graph are adjacent (by the definition of the edges of this graph).
Qed.}. Hence, the vertices $\gamma_{0},\gamma_{1},\ldots,\gamma_{L}$ form a
path. This path connects $\sigma$ to $\tau$ (since it begins at $\gamma
_{0}=\sigma$ and ends at $\gamma_{L}=\tau$), and has length $L$. Thus, there
exists a path between $\sigma$ and $\tau$ which has length $L$ (namely, the
path formed by the vertices $\gamma_{0},\gamma_{1},\ldots,\gamma_{L}$).

We shall now show that $L$ is the smallest length of a path between $\sigma$
and $\tau$. Indeed, let $\mathbf{d}$ be any path between $\sigma$ and $\tau$.
We shall show that the length of $\mathbf{d}$ is $\geq L$.

The path $\mathbf{d}$ is a path between $\sigma$ and $\tau$. Hence, we can
write the path $\mathbf{d}$ in the form $\mathbf{d}=\left(  \delta_{0}%
,\delta_{1},\ldots,\delta_{M}\right)  $ for some $\delta_{0},\delta_{1}%
,\ldots,\delta_{M}\in S_{n}$ with $\delta_{0}=\sigma$ and $\delta_{M}=\tau$.
Consider these $\delta_{0},\delta_{1},\ldots,\delta_{M}\in S_{n}$.

For every $i\in\left\{  1,2,\ldots,M\right\}  $, there exists a $k\in\left\{
1,2,\ldots,n-1\right\}  $ such that $\delta_{i}=\delta_{i-1}\circ s_{k}%
$\ \ \ \ \footnote{\textit{Proof.} Let $i\in\left\{  1,2,\ldots,M\right\}  $.
Then, the vertices $\delta_{i}$ and $\delta_{i-1}$ of the (undirected) $n$-th
right Bruhat graph are adjacent (because they are two consecutive vertices on
the path $\left(  \delta_{0},\delta_{1},\ldots,\delta_{M}\right)  =\mathbf{d}%
$). In other words, there exists a $k\in\left\{  1,2,\ldots,n-1\right\}  $
such that $\delta_{i}=\delta_{i-1}\circ s_{k}$ (by the definition of the edges
of this graph). Qed.}. We denote this $k$ by $k_{i}$. Thus, for every
$i\in\left\{  1,2,\ldots,M\right\}  $, we have defined a $k_{i}\in\left\{
1,2,\ldots,n-1\right\}  $ such that $\delta_{i}=\delta_{i-1}\circ s_{k_{i}}$.

Now, every $j\in\left\{  0,1,\ldots,M\right\}  $ satisfies%
\begin{equation}
\delta_{j}=\sigma\circ\left(  s_{k_{1}}\circ s_{k_{2}}\circ\cdots\circ
s_{k_{j}}\right)  \label{sol.ps4.0.5}%
\end{equation}
\footnote{\textit{Proof of (\ref{sol.ps4.0.5}):} We shall prove
(\ref{sol.ps4.0.5}) by induction over $j$:
\par
\textit{Induction base:} We have $\delta_{0}=\sigma$. Compared with
$\sigma\circ\underbrace{\left(  s_{k_{1}}\circ s_{k_{2}}\circ\cdots\circ
s_{k_{0}}\right)  }_{=\left(  \text{a composition of }0\text{ permutations}%
\right)  =\operatorname*{id}}=\sigma$, this yields $\delta_{0}=\sigma
\circ\left(  s_{k_{1}}\circ s_{k_{2}}\circ\cdots\circ s_{k_{0}}\right)  $. In
other words, (\ref{sol.ps4.0.5}) holds for $j=0$. This completes the induction
base.
\par
\textit{Induction step:} Let $J\in\left\{  0,1,\ldots,M\right\}  $ be
positive. Assume that (\ref{sol.ps4.0.5}) holds for $j=J-1$. We need to show
that (\ref{sol.ps4.0.5}) holds for $j=J$.
\par
We have $\delta_{J}=\sigma\circ\left(  s_{k_{1}}\circ s_{k_{2}}\circ
\cdots\circ s_{k_{J-1}}\right)  $ (since (\ref{sol.ps4.0.5}) holds for
$j=J-1$). Now, recall that $\delta_{i}=\delta_{i-1}\circ s_{k_{i}}$ for every
$i\in\left\{  1,2,\ldots,M\right\}  $. Applying this to $i=J$, we obtain%
\begin{align*}
\delta_{J}  &  =\underbrace{\delta_{J-1}}_{=\sigma\circ\left(  s_{k_{1}}\circ
s_{k_{2}}\circ\cdots\circ s_{k_{J-1}}\right)  }\circ s_{k_{J}}%
\ \ \ \ \ \ \ \ \ \ \left(  \text{since }J\in\left\{  1,2,\ldots,M\right\}
\text{ (since }J\in\left\{  0,1,\ldots,M\right\}  \text{ is positive)}\right)
\\
&  =\sigma\circ\underbrace{\left(  s_{k_{1}}\circ s_{k_{2}}\circ\cdots\circ
s_{k_{J-1}}\right)  \circ s_{k_{J}}}_{=s_{k_{1}}\circ s_{k_{2}}\circ
\cdots\circ s_{k_{J}}}=\sigma\circ\left(  s_{k_{1}}\circ s_{k_{2}}\circ
\cdots\circ s_{k_{J}}\right)  .
\end{align*}
In other words, (\ref{sol.ps4.0.5}) holds for $j=J$. This completes the
induction step. Thus, (\ref{sol.ps4.0.5}) is proven by induction.}. Applying
this to $j=M$, we obtain%
\[
\delta_{M}=\sigma\circ\left(  s_{k_{1}}\circ s_{k_{2}}\circ\cdots\circ
s_{k_{M}}\right)  .
\]
Compared with $\delta_{M}=\tau$, this yields%
\[
\tau=\sigma\circ\left(  s_{k_{1}}\circ s_{k_{2}}\circ\cdots\circ s_{k_{M}%
}\right)  ,
\]
so that
\begin{align*}
\sigma^{-1}\circ\underbrace{\tau}_{=\sigma\circ\left(  s_{k_{1}}\circ
s_{k_{2}}\circ\cdots\circ s_{k_{M}}\right)  }  &  =\underbrace{\sigma
^{-1}\circ\sigma}_{=\operatorname*{id}}\circ\left(  s_{k_{1}}\circ s_{k_{2}%
}\circ\cdots\circ s_{k_{M}}\right) \\
&  =s_{k_{1}}\circ s_{k_{2}}\circ\cdots\circ s_{k_{M}}.
\end{align*}
Therefore, the permutation $\sigma^{-1}\circ\tau$ can be written as a
composition of $M$ permutations of the form $s_{k}$ (with $k\in\left\{
1,2,\ldots,n-1\right\}  $) (namely, of the $M$ permutations $s_{k_{1}}$,
$s_{k_{2}}$, $\ldots$, $s_{k_{M}}$).

Now, we recall that $\ell\left(  \sigma^{-1}\circ\tau\right)  $ is the
smallest $N\in\mathbb{N}$ such that $\sigma^{-1}\circ\tau$ can be written as a
composition of $N$ permutations of the form $s_{k}$ (with $k\in\left\{
1,2,\ldots,n-1\right\}  $). Hence, if $N\in\mathbb{N}$ is such that
$\sigma^{-1}\circ\tau$ can be written as a composition of $N$ permutations of
the form $s_{k}$ (with $k\in\left\{  1,2,\ldots,n-1\right\}  $), then
$N\geq\ell\left(  \sigma^{-1}\circ\tau\right)  $. We can apply this to $N=M$
(because $\sigma^{-1}\circ\tau$ can be written as a composition of $M$
permutations of the form $s_{k}$ (with $k\in\left\{  1,2,\ldots,n-1\right\}
$)), and thus obtain $M\geq\ell\left(  \sigma^{-1}\circ\tau\right)  =L$.

But the length of the path $\mathbf{d}$ is $M$ (since $\mathbf{d}=\left(
\delta_{0},\delta_{1},\ldots,\delta_{M}\right)  $). Hence, the length of the
path $\mathbf{d}$ is $\geq L$ (since $M\geq L$).

Let us now forget that we fixed $\mathbf{d}$. We thus have shown that if
$\mathbf{d}$ is any path between $\sigma$ and $\tau$, then the length of the
path $\mathbf{d}$ is $\geq L$. In other words, every path between $\sigma$ and
$\tau$ has length $\geq L$.

Altogether, we have proven the following two statements:

\begin{itemize}
\item There exists a path of length $L$ between $\sigma$ and $\tau$.

\item Every path between $\sigma$ and $\tau$ has length $\geq L$.
\end{itemize}

Therefore, $L$ is the smallest length of a path between $\sigma$ and $\tau$.
In other words, $\ell\left(  \sigma^{-1}\circ\tau\right)  $ is the smallest
length of a path between $\sigma$ and $\tau$ (since $L=\ell\left(  \sigma
^{-1}\circ\tau\right)  $). Exercise \ref{exe.ps4.0} is solved.
\end{proof}
\end{verlong}

\subsection{\label{sect.sols.transpose.code}Solution to Exercise
\ref{exe.transpos.code}}

\subsubsection{Preparations}

Before we solve Exercise \ref{exe.transpos.code}, let us prepare with some
simple lemmas and notations. The following notation will be used throughout
Section \ref{sect.sols.transpose.code}:

\begin{definition}
\label{def.set12...m}Whenever $m$ is an integer, we shall use the notation
$\left[  m\right]  $ for the set $\left\{  1,2,\ldots,m\right\}  $. (This is
an empty set when $m\leq0$.)
\end{definition}

Notice that if $a$ and $b$ are two integers satisfying $a\leq b$, then
$\left[  a\right]  \subseteq\left[  b\right]  $.

The following is a simple property of the permutations $t_{i,j}$ defined in
Definition \ref{def.transpos.ii}:

\begin{lemma}
\label{lem.sol.transpose.code.tij1}Let $n\in\mathbb{N}$. Let $i$ and $j$ be
two elements of $\left[  n\right]  $.

\textbf{(a)} We have $t_{i,j}\left(  i\right)  =j$.

\textbf{(b)} We have $t_{i,j}\left(  j\right)  =i$.

\textbf{(c)} We have $t_{i,j}\left(  k\right)  =k$ for each $k\in\left[
n\right]  \setminus\left\{  i,j\right\}  $.

\textbf{(d)} We have $t_{i,j}\circ t_{i,j}=\operatorname*{id}$.
\end{lemma}

\begin{vershort}
\begin{proof}
[Proof of Lemma \ref{lem.sol.transpose.code.tij1}.]Lemma
\ref{lem.sol.transpose.code.tij1} follows from the definition of $t_{i,j}$
(given in Definition \ref{def.transpos} and in Definition
\ref{def.transpos.ii}).
\end{proof}
\end{vershort}

\begin{verlong}
\begin{proof}
[Proof of Lemma \ref{lem.sol.transpose.code.tij1}.]We have $\left[  n\right]
=\left\{  1,2,\ldots,n\right\}  $ (by the definition of $\left[  n\right]  $).
We know that $i$ and $j$ are two elements of $\left[  n\right]  $. In other
words, $i$ and $j$ are two elements of $\left\{  1,2,\ldots,n\right\}  $
(since $\left[  n\right]  =\left\{  1,2,\ldots,n\right\}  $).

\textbf{(a)} We are in one of the following two cases:

\textit{Case 1:} We have $i=j$.

\textit{Case 2:} We have $i\neq j$.

Let us first consider Case 1. In this case, we have $i=j$. Hence, Definition
\ref{def.transpos.ii} shows that $t_{i,j}=\operatorname*{id}$. Thus,
$\underbrace{t_{i,j}}_{=\operatorname*{id}}\left(  i\right)
=\operatorname*{id}\left(  i\right)  =i=j$. Hence, Lemma
\ref{lem.sol.transpose.code.tij1} \textbf{(a)} is proven in Case 1.

Let us now consider Case 2. In this case, we have $i\neq j$. In other words,
the elements $i$ and $j$ are distinct. Thus, Definition \ref{def.transpos}
shows that $t_{i,j}$ is the permutation in $S_{n}$ which swaps $i$ with $j$
while leaving all other elements of $\left\{  1,2,\ldots,n\right\}  $
unchanged. Hence, the permutation $t_{i,j}$ swaps $i$ with $j$. In other
words, it satisfies $t_{i,j}\left(  i\right)  =j$ and $t_{i,j}\left(
j\right)  =i$. In particular, we have $t_{i,j}\left(  i\right)  =j$. Hence,
Lemma \ref{lem.sol.transpose.code.tij1} \textbf{(a)} is proven in Case 2.

We have now proven Lemma \ref{lem.sol.transpose.code.tij1} \textbf{(a)} in
each of the two Cases 1 and 2. Since these two Cases cover all possibilities,
we thus conclude that Lemma \ref{lem.sol.transpose.code.tij1} \textbf{(a)}
always holds. Thus, Lemma \ref{lem.sol.transpose.code.tij1} \textbf{(a)} is proven.

\textbf{(b)} We are in one of the following two cases:

\textit{Case 1:} We have $i=j$.

\textit{Case 2:} We have $i\neq j$.

Let us first consider Case 1. In this case, we have $i=j$. Hence, Definition
\ref{def.transpos.ii} shows that $t_{i,j}=\operatorname*{id}$. Thus,
$\underbrace{t_{i,j}}_{=\operatorname*{id}}\left(  j\right)
=\operatorname*{id}\left(  j\right)  =j=i$. Hence, Lemma
\ref{lem.sol.transpose.code.tij1} \textbf{(b)} is proven in Case 1.

Let us now consider Case 2. In this case, we have $i\neq j$. In other words,
the elements $i$ and $j$ are distinct. Thus, Definition \ref{def.transpos}
shows that $t_{i,j}$ is the permutation in $S_{n}$ which swaps $i$ with $j$
while leaving all other elements of $\left\{  1,2,\ldots,n\right\}  $
unchanged. Hence, the permutation $t_{i,j}$ swaps $i$ with $j$. In other
words, it satisfies $t_{i,j}\left(  i\right)  =j$ and $t_{i,j}\left(
j\right)  =i$. In particular, we have $t_{i,j}\left(  j\right)  =i$. Hence,
Lemma \ref{lem.sol.transpose.code.tij1} \textbf{(b)} is proven in Case 2.

We have now proven Lemma \ref{lem.sol.transpose.code.tij1} \textbf{(b)} in
each of the two Cases 1 and 2. Since these two Cases cover all possibilities,
we thus conclude that Lemma \ref{lem.sol.transpose.code.tij1} \textbf{(b)}
always holds. Thus, Lemma \ref{lem.sol.transpose.code.tij1} \textbf{(b)} is proven.

\textbf{(c)} Let $k\in\left[  n\right]  \setminus\left\{  i,j\right\}  $. We
must prove that $t_{i,j}\left(  k\right)  =k$.

We have $k\in\underbrace{\left[  n\right]  }_{=\left\{  1,2,\ldots,n\right\}
}\setminus\left\{  i,j\right\}  =\left\{  1,2,\ldots,n\right\}  \setminus
\left\{  i,j\right\}  $. In other words, $k$ is an element of $\left\{
1,2,\ldots,n\right\}  $ other than $i$ and $j$.

We are in one of the following two cases:

\textit{Case 1:} We have $i=j$.

\textit{Case 2:} We have $i\neq j$.

Let us first consider Case 1. In this case, we have $i=j$. Hence, Definition
\ref{def.transpos.ii} shows that $t_{i,j}=\operatorname*{id}$. Thus,
$\underbrace{t_{i,j}}_{=\operatorname*{id}}\left(  k\right)
=\operatorname*{id}\left(  k\right)  =k$. Hence, Lemma
\ref{lem.sol.transpose.code.tij1} \textbf{(c)} is proven in Case 1.

Let us now consider Case 2. In this case, we have $i\neq j$. In other words,
the elements $i$ and $j$ are distinct. Thus, Definition \ref{def.transpos}
shows that $t_{i,j}$ is the permutation in $S_{n}$ which swaps $i$ with $j$
while leaving all other elements of $\left\{  1,2,\ldots,n\right\}  $
unchanged. Hence, the permutation $t_{i,j}$ leaves all elements of $\left\{
1,2,\ldots,n\right\}  $ other than $i$ and $j$ unchanged. In other words, it
satisfies $t_{i,j}\left(  p\right)  =p$ whenever $p$ is an element of
$\left\{  1,2,\ldots,n\right\}  $ other than $i$ and $j$. We can apply this to
$p=k$ (since $k$ is an element of $\left\{  1,2,\ldots,n\right\}  $ other than
$i$ and $j$), and thus obtain $t_{i,j}\left(  k\right)  =k$. Hence, Lemma
\ref{lem.sol.transpose.code.tij1} \textbf{(c)} is proven in Case 2.

We have now proven Lemma \ref{lem.sol.transpose.code.tij1} \textbf{(c)} in
each of the two Cases 1 and 2. Since these two Cases cover all possibilities,
we thus conclude that Lemma \ref{lem.sol.transpose.code.tij1} \textbf{(c)}
always holds. Thus, Lemma \ref{lem.sol.transpose.code.tij1} \textbf{(c)} is proven.

\textbf{(d)} Clearly, $t_{i,j}$ is a permutation in $S_{n}$, thus a bijection
$\left\{  1,2,\ldots,n\right\}  \rightarrow\left\{  1,2,\ldots,n\right\}  $.
Hence, $t_{i,j}\circ t_{i,j}$ is a map $\left\{  1,2,\ldots,n\right\}
\rightarrow\left\{  1,2,\ldots,n\right\}  $ as well.

Let $k\in\left\{  1,2,\ldots,n\right\}  $. We shall show that $\left(
t_{i,j}\circ t_{i,j}\right)  \left(  k\right)  =k$.

We have $k\in\left\{  1,2,\ldots,n\right\}  =\left[  n\right]  $. We are in
one of the following three cases:

\textit{Case 1:} We have $k=i$.

\textit{Case 2:} We have $k=j$.

\textit{Case 3:} We have neither $k=i$ nor $k=j$.

Let us first consider Case 1. In this case, we have $k=i$. Thus,%
\begin{align*}
\left(  t_{i,j}\circ t_{i,j}\right)  \left(  \underbrace{k}_{=i}\right)   &
=\left(  t_{i,j}\circ t_{i,j}\right)  \left(  i\right)  =t_{i,j}\left(
\underbrace{t_{i,j}\left(  i\right)  }_{\substack{=j\\\text{(by Lemma
\ref{lem.sol.transpose.code.tij1} \textbf{(a)})}}}\right)  =t_{i,j}\left(
j\right) \\
&  =i\ \ \ \ \ \ \ \ \ \ \left(  \text{by Lemma
\ref{lem.sol.transpose.code.tij1} \textbf{(b)}}\right) \\
&  =k.
\end{align*}
Thus, $\left(  t_{i,j}\circ t_{i,j}\right)  \left(  k\right)  =k$ is proven in
Case 1.

Let us now consider Case 2. In this case, we have $k=j$. Thus,%
\begin{align*}
\left(  t_{i,j}\circ t_{i,j}\right)  \left(  \underbrace{k}_{=j}\right)   &
=\left(  t_{i,j}\circ t_{i,j}\right)  \left(  j\right)  =t_{i,j}\left(
\underbrace{t_{i,j}\left(  j\right)  }_{\substack{=i\\\text{(by Lemma
\ref{lem.sol.transpose.code.tij1} \textbf{(b)})}}}\right)  =t_{i,j}\left(
i\right) \\
&  =j\ \ \ \ \ \ \ \ \ \ \left(  \text{by Lemma
\ref{lem.sol.transpose.code.tij1} \textbf{(a)}}\right) \\
&  =k.
\end{align*}
Thus, $\left(  t_{i,j}\circ t_{i,j}\right)  \left(  k\right)  =k$ is proven in
Case 2.

Let us finally consider Case 3. In this case, we have neither $k=i$ nor $k=j$.
Hence, we have $k\notin\left\{  i,j\right\}  $%
\ \ \ \ \footnote{\textit{Proof.} Assume the contrary. Thus, $k\in\left\{
i,j\right\}  $. In other words, $k$ is either $i$ or $j$. In other words, we
have either $k=i$ or $k=j$. This contradicts the fact that we have neither
$k=i$ nor $k=j$. This contradiction shows that our assumption was wrong.
Qed.}. Combining $k\in\left[  n\right]  $ with $k\notin\left\{  i,j\right\}
$, we obtain $k\in\left[  n\right]  \setminus\left\{  i,j\right\}  $. Hence,
$t_{i,j}\left(  k\right)  =k$ (by Lemma \ref{lem.sol.transpose.code.tij1}
\textbf{(c)}). Now,
\[
\left(  t_{i,j}\circ t_{i,j}\right)  \left(  k\right)  =t_{i,j}\left(
\underbrace{t_{i,j}\left(  k\right)  }_{=k}\right)  =t_{i,j}\left(  k\right)
=k.
\]
Thus, $\left(  t_{i,j}\circ t_{i,j}\right)  \left(  k\right)  =k$ is proven in
Case 3.

We now have proven $\left(  t_{i,j}\circ t_{i,j}\right)  \left(  k\right)  =k$
in each of the three Cases 1, 2 and 3. Since these three Cases cover all
possibilities, we thus conclude that $\left(  t_{i,j}\circ t_{i,j}\right)
\left(  k\right)  =k$ always holds. Thus, $\left(  t_{i,j}\circ t_{i,j}%
\right)  \left(  k\right)  =k=\operatorname*{id}\left(  k\right)  $.

Now, forget that we fixed $k$. We thus have proven that $\left(  t_{i,j}\circ
t_{i,j}\right)  \left(  k\right)  =\operatorname*{id}\left(  k\right)  $ for
each $k\in\left\{  1,2,\ldots,n\right\}  $. Hence, $t_{i,j}\circ
t_{i,j}=\operatorname*{id}$ (since both $t_{i,j}\circ t_{i,j}$ and
$\operatorname*{id}$ are maps $\left\{  1,2,\ldots,n\right\}  \rightarrow
\left\{  1,2,\ldots,n\right\}  $). This proves Lemma
\ref{lem.sol.transpose.code.tij1} \textbf{(d)}.
\end{proof}
\end{verlong}

\begin{lemma}
\label{lem.sol.transpose.code.tij2}Let $n\in\mathbb{N}$. Let $k\in\left[
n\right]  $. Let $\sigma\in S_{n}$ be such that%
\begin{equation}
\left(  \sigma\left(  i\right)  =i\text{ for each }i\in\left\{  k+1,k+2,\ldots
,n\right\}  \right)  . \label{eq.lem.sol.transpose.code.tij2.ass}%
\end{equation}
Let $g=\sigma^{-1}\left(  k\right)  $. Then:

\textbf{(a)} We have $g\in\left[  k\right]  $.

\textbf{(b)} We have $\left(  \sigma\circ t_{k,g}\right)  \left(  i\right)
=i$ for each $i\in\left\{  k,k+1,\ldots,n\right\}  $.
\end{lemma}

\begin{vershort}
\begin{proof}
[Proof of Lemma \ref{lem.sol.transpose.code.tij2}.]We have $\sigma\in S_{n}$.
In other words, $\sigma$ is a permutation of $\left\{  1,2,\ldots,n\right\}  $
(since $S_{n}$ is the set of all permutations of $\left\{  1,2,\ldots
,n\right\}  $), that is, a bijection $\left[  n\right]  \rightarrow\left[
n\right]  $. Hence, the inverse $\sigma^{-1}$ of $\sigma$ is well-defined.

From $g=\sigma^{-1}\left(  k\right)  $, we obtain $\sigma\left(  g\right)  =k$.

\textbf{(a)} Assume the contrary. Thus, $g\notin\left[  k\right]  $. Combining
$g\in\left[  n\right]  $ with $g\notin\left[  k\right]  $, we obtain
$g\in\left[  n\right]  \setminus\left[  k\right]  =\left\{  k+1,k+2,\ldots
,n\right\}  $. Hence, (\ref{eq.lem.sol.transpose.code.tij2.ass}) (applied to
$i=g$) yields $\sigma\left(  g\right)  =g$. Hence, $g=\sigma\left(  g\right)
=k\in\left[  k\right]  $. This contradicts $\sigma\left(  k\right)
\notin\left[  k\right]  $. This contradiction shows that our assumption was
wrong. Hence, Lemma \ref{lem.sol.transpose.code.tij2} \textbf{(a)} is proven.

\textbf{(b)} Let $i\in\left\{  k,k+1,\ldots,n\right\}  $. We must prove that
$\left(  \sigma\circ t_{k,g}\right)  \left(  i\right)  =i$.

We are in one of the following two cases:

\textit{Case 1:} We have $i\neq k$.

\textit{Case 2:} We have $i=k$.

Let us first consider Case 1. In this case, we have $i\neq k$. Combining
$i\in\left\{  k,k+1,\ldots,n\right\}  $ with $i\neq k$, we obtain
$i\in\left\{  k,k+1,\ldots,n\right\}  \setminus\left\{  k\right\}  =\left\{
k+1,k+2,\ldots,n\right\}  $. Hence, (\ref{eq.lem.sol.transpose.code.tij2.ass})
shows that $\sigma\left(  i\right)  =i$.

But $i\notin\left\{  k,g\right\}  $\ \ \ \ \footnote{\textit{Proof.} Assume
the contrary. Thus, $i\in\left\{  k,g\right\}  $. Combining this with $i\neq
k$, we obtain $i\in\left\{  k,g\right\}  \setminus\left\{  k\right\}
\subseteq\left\{  g\right\}  $. Thus, $i=g\in\left[  k\right]  $ (by Lemma
\ref{lem.sol.transpose.code.tij2} \textbf{(a)}). Hence, $i\leq k$. But
$i\in\left\{  k+1,k+2,\ldots,n\right\}  $, so that $i\geq k+1>k$. This
contradicts $i\leq k$. This contradiction shows that our assumption was wrong.
Qed.}. Combining $i\in\left\{  k,k+1,\ldots,n\right\}  \subseteq\left[
n\right]  $ with $i\notin\left\{  k,g\right\}  $, we obtain $i\in\left[
n\right]  \setminus\left\{  k,g\right\}  $. Hence, Lemma
\ref{lem.sol.transpose.code.tij1} \textbf{(c)} (applied to $k$, $g$ and $i$
instead of $i$, $j$ and $k$) shows that $t_{k,g}\left(  i\right)  =i$.

Now, $\left(  \sigma\circ t_{k,g}\right)  \left(  i\right)  =\sigma\left(
\underbrace{t_{k,g}\left(  i\right)  }_{=i}\right)  =\sigma\left(  i\right)
=i$. Thus, $\left(  \sigma\circ t_{k,g}\right)  \left(  i\right)  =i$ is
proven in Case 1.

Let us now consider Case 2. In this case, we have $i=k$. Lemma
\ref{lem.sol.transpose.code.tij1} \textbf{(a)} (applied to $k$ and $g$ instead
of $i$ and $j$) yields $t_{k,g}\left(  k\right)  =g$. Thus, $\left(
\sigma\circ t_{k,g}\right)  \left(  \underbrace{i}_{=k}\right)  =\left(
\sigma\circ t_{k,g}\right)  \left(  k\right)  =\sigma\left(
\underbrace{t_{k,g}\left(  k\right)  }_{=g}\right)  =\sigma\left(  g\right)
=k=i$. Thus, $\left(  \sigma\circ t_{k,g}\right)  \left(  i\right)  =i$ is
proven in Case 2.

We now have proven $\left(  \sigma\circ t_{k,g}\right)  \left(  i\right)  =i$
in each of the two Cases 1 and 2. Hence, $\left(  \sigma\circ t_{k,g}\right)
\left(  i\right)  =i$ always holds. Thus, Lemma
\ref{lem.sol.transpose.code.tij2} \textbf{(b)} is proven.
\end{proof}
\end{vershort}

\begin{verlong}
\begin{proof}
[Proof of Lemma \ref{lem.sol.transpose.code.tij2}.]The definition of $\left[
k\right]  $ yields $\left[  k\right]  =\left\{  1,2,\ldots,k\right\}  $.

We have $k\in\left[  n\right]  =\left\{  1,2,\ldots,n\right\}  $ (by the
definition of $\left[  n\right]  $). Hence, $k$ is a positive integer. Thus,
$k\in\left\{  1,2,\ldots,k\right\}  =\left[  k\right]  $.

We have $\sigma\in S_{n}$. In other words, $\sigma$ is a permutation of
$\left\{  1,2,\ldots,n\right\}  $ (since $S_{n}$ is the set of all
permutations of $\left\{  1,2,\ldots,n\right\}  $). In other words, $\sigma$
is a permutation of $\left[  n\right]  $ (since $\left[  n\right]  =\left\{
1,2,\ldots,n\right\}  $). In other words, $\sigma$ is a bijection $\left[
n\right]  \rightarrow\left[  n\right]  $. Hence, this map $\sigma$ is
invertible, and its inverse $\sigma^{-1}$ is also a bijection $\left[
n\right]  \rightarrow\left[  n\right]  $.

We have $g=\sigma^{-1}\left(  k\right)  \in\left[  n\right]  $ (since
$\sigma^{-1}$ is a map $\left[  n\right]  \rightarrow\left[  n\right]  $).
Also, $\sigma\left(  \underbrace{g}_{=\sigma^{-1}\left(  k\right)  }\right)
=\sigma\left(  \sigma^{-1}\left(  k\right)  \right)  =k$.

From $k\in\left\{  1,2,\ldots,n\right\}  $, we obtain $k\geq1$. Hence,
$\left\{  k,k+1,\ldots,n\right\}  \subseteq\left\{  1,2,\ldots,n\right\}
=\left[  n\right]  $.

\textbf{(a)} Assume the contrary. Thus, $g\notin\left[  k\right]  $.

Combining $g\in\left[  n\right]  =\left\{  1,2,\ldots,n\right\}  $ with
$g\notin\left[  k\right]  =\left\{  1,2,\ldots,k\right\}  $, we obtain
\[
g\in\left\{  1,2,\ldots,n\right\}  \setminus\left\{  1,2,\ldots,k\right\}
=\left\{  k+1,k+2,\ldots,n\right\}  .
\]
Hence, (\ref{eq.lem.sol.transpose.code.tij2.ass}) (applied to $i=g$) yields
$\sigma\left(  g\right)  =g$. Hence, $g=\sigma\left(  g\right)  =k\in\left[
k\right]  $. This contradicts $\sigma\left(  k\right)  \notin\left[  k\right]
$.

This contradiction shows that our assumption was wrong. Hence, Lemma
\ref{lem.sol.transpose.code.tij2} \textbf{(a)} is proven.

\textbf{(b)} Let $i\in\left\{  k,k+1,\ldots,n\right\}  $. We must prove that
$\left(  \sigma\circ t_{k,g}\right)  \left(  i\right)  =i$.

We are in one of the following two cases:

\textit{Case 1:} We have $i\neq k$.

\textit{Case 2:} We have $i=k$.

Let us first consider Case 1. In this case, we have $i\neq k$. Combining
$i\in\left\{  k,k+1,\ldots,n\right\}  $ with $i\neq k$, we obtain
$i\in\left\{  k,k+1,\ldots,n\right\}  \setminus\left\{  k\right\}  =\left\{
k+1,k+2,\ldots,n\right\}  $. Hence, (\ref{eq.lem.sol.transpose.code.tij2.ass})
shows that $\sigma\left(  i\right)  =i$.

But $i\notin\left\{  k,g\right\}  $\ \ \ \ \footnote{\textit{Proof.} Assume
the contrary. Thus, $i\in\left\{  k,g\right\}  $. Combining this with $i\neq
k$, we obtain $i\in\left\{  k,g\right\}  \setminus\left\{  k\right\}
\subseteq\left\{  g\right\}  $. Thus, $i=g\in\left[  k\right]  $ (by Lemma
\ref{lem.sol.transpose.code.tij2} \textbf{(a)}). Hence, $i\in\left[  k\right]
=\left\{  1,2,\ldots,k\right\}  $, so that $i\leq k$. But $i\in\left\{
k+1,k+2,\ldots,n\right\}  $, so that $i\geq k+1>k$. This contradicts $i\leq
k$. This contradiction shows that our assumption was wrong. Qed.}. Combining
$i\in\left\{  k,k+1,\ldots,n\right\}  \subseteq\left[  n\right]  $ with
$i\notin\left\{  k,g\right\}  $, we obtain $i\in\left[  n\right]
\setminus\left\{  k,g\right\}  $. Hence, Lemma
\ref{lem.sol.transpose.code.tij1} \textbf{(c)} (applied to $k$, $g$ and $i$
instead of $i$, $j$ and $k$) shows that $t_{k,g}\left(  i\right)  =i$.

Now, $\left(  \sigma\circ t_{k,g}\right)  \left(  i\right)  =\sigma\left(
\underbrace{t_{k,g}\left(  i\right)  }_{=i}\right)  =\sigma\left(  i\right)
=i$. Thus, $\left(  \sigma\circ t_{k,g}\right)  \left(  i\right)  =i$ is
proven in Case 1.

Let us now consider Case 2. In this case, we have $i=k$. On the other hand,
Lemma \ref{lem.sol.transpose.code.tij1} \textbf{(a)} (applied to $k$ and $g$
instead of $i$ and $j$) yields $t_{k,g}\left(  k\right)  =g$. Thus, $\left(
\sigma\circ t_{k,g}\right)  \left(  \underbrace{i}_{=k}\right)  =\left(
\sigma\circ t_{k,g}\right)  \left(  k\right)  =\sigma\left(
\underbrace{t_{k,g}\left(  k\right)  }_{=g}\right)  =\sigma\left(  g\right)
=k=i$. Thus, $\left(  \sigma\circ t_{k,g}\right)  \left(  i\right)  =i$ is
proven in Case 2.

We now have proven $\left(  \sigma\circ t_{k,g}\right)  \left(  i\right)  =i$
in each of the two Cases 1 and 2. Since these two Cases cover all
possibilities, we thus conclude that $\left(  \sigma\circ t_{k,g}\right)
\left(  i\right)  =i$ always holds. Thus, Lemma
\ref{lem.sol.transpose.code.tij2} \textbf{(b)} is proven.
\end{proof}
\end{verlong}

\begin{noncompile}
The following material has been commented out since it is not used.

\begin{lemma}
\label{lem.sol.transpose.code.tij3}Let $n\in\mathbb{N}$. Let $k\in\left[
n\right]  $. Let $\sigma\in S_{n}$ be such that%
\begin{equation}
\left(  \sigma\left(  i\right)  =i\text{ for each }i\in\left\{  k+1,k+2,\ldots
,n\right\}  \right)  . \label{eq.lem.sol.transpose.code.tij3.ass}%
\end{equation}
Then:

\textbf{(a)} We have $\sigma\left(  k\right)  \in\left[  k\right]  $.

\textbf{(b)} We have $\left(  t_{k,\sigma\left(  k\right)  }\circ
\sigma\right)  \left(  i\right)  =i$ for each $i\in\left\{  k,k+1,\ldots
,n\right\}  $.

\textbf{(c)} If $u\in\left[  n\right]  $ is such that $\left(  t_{k,u}%
\circ\sigma\right)  \left(  i\right)  =i$ for each $i\in\left\{
k,k+1,\ldots,n\right\}  $, then $u=\sigma\left(  k\right)  $.
\end{lemma}

VERSHORT PROOF:

\begin{proof}
[Proof of Lemma \ref{lem.sol.transpose.code.tij3}.]We have $\sigma\in S_{n}$.
In other words, $\sigma$ is a permutation of $\left[  n\right]  $ (since
$S_{n}$ is the set of all permutations of $\left[  n\right]  $). Thus, the map
$\sigma$ is bijective, hence injective.

\textbf{(a)} Assume the contrary. Thus, $\sigma\left(  k\right)  \notin\left[
k\right]  $.

Combining $\sigma\left(  k\right)  \in\left[  n\right]  $ with $\sigma\left(
k\right)  \notin\left[  k\right]  $, we obtain $\sigma\left(  k\right)
\in\left[  n\right]  \setminus\left[  k\right]  =\left\{  k+1,k+2,\ldots
,n\right\}  $. Hence, (\ref{eq.lem.sol.transpose.code.tij3.ass}) (applied to
$i=\sigma\left(  k\right)  $) yields $\sigma\left(  \sigma\left(  k\right)
\right)  =\sigma\left(  k\right)  $. Since $\sigma$ is injective, we thus
conclude that $\sigma\left(  k\right)  =k\in\left[  k\right]  $. This
contradicts $\sigma\left(  k\right)  \notin\left[  k\right]  $.

This contradiction shows that our assumption was wrong. Hence, Lemma
\ref{lem.sol.transpose.code.tij3} \textbf{(a)} is proven.

\textbf{(b)} Let $i\in\left\{  k,k+1,\ldots,n\right\}  $. We must prove that
$\left(  t_{k,\sigma\left(  k\right)  }\circ\sigma\right)  \left(  i\right)
=i$.

We are in one of the following two cases:

\textit{Case 1:} We have $i\neq k$.

\textit{Case 2:} We have $i=k$.

Let us first consider Case 1. In this case, we have $i\neq k$. Combining
$i\in\left\{  k,k+1,\ldots,n\right\}  $ with $i\neq k$, we obtain
$i\in\left\{  k,k+1,\ldots,n\right\}  \setminus\left\{  k\right\}  =\left\{
k+1,k+2,\ldots,n\right\}  $. Hence, (\ref{eq.lem.sol.transpose.code.tij3.ass})
shows that $\sigma\left(  i\right)  =i$.

But $i\notin\left\{  k,\sigma\left(  k\right)  \right\}  $%
\ \ \ \ \footnote{\textit{Proof.} Assume the contrary. Thus, $i\in\left\{
k,\sigma\left(  k\right)  \right\}  $. Combining this with $i\neq k$, we
obtain $i\in\left\{  k,\sigma\left(  k\right)  \right\}  \setminus\left\{
k\right\}  \subseteq\left\{  \sigma\left(  k\right)  \right\}  $. Thus,
$i=\sigma\left(  k\right)  \in\left[  k\right]  $ (by Lemma
\ref{lem.sol.transpose.code.tij3} \textbf{(a)}). Hence, $i\in\left[  k\right]
=\left\{  1,2,\ldots,k\right\}  $, so that $i\leq k$. But $i\in\left\{
k+1,k+2,\ldots,n\right\}  $, so that $i\geq k+1>k$. This contradicts $i\leq
k$. This contradiction shows that our assumption was wrong. Qed.}. Combining
$i\in\left\{  k,k+1,\ldots,n\right\}  \subseteq\left\{  1,2,\ldots,n\right\}
=\left[  n\right]  $ with $i\notin\left\{  k,\sigma\left(  k\right)  \right\}
$, we obtain $i\in\left[  n\right]  \setminus\left\{  k,\sigma\left(
k\right)  \right\}  $. Hence, Lemma \ref{lem.sol.transpose.code.tij1}
\textbf{(c)} (applied to $k$, $\sigma\left(  k\right)  $ and $i$ instead of
$i$, $j$ and $k$) shows that $t_{k,\sigma\left(  k\right)  }\left(  i\right)
=i$.

Now, $\left(  t_{k,\sigma\left(  k\right)  }\circ\sigma\right)  \left(
i\right)  =t_{k,\sigma\left(  k\right)  }\left(  \underbrace{\sigma\left(
i\right)  }_{=i}\right)  =t_{k,\sigma\left(  k\right)  }\left(  i\right)  =i$.
Thus, $\left(  t_{k,\sigma\left(  k\right)  }\circ\sigma\right)  \left(
i\right)  =i$ is proven in Case 1.

Let us now consider Case 2. In this case, we have $i=k$. On the other hand,
Lemma \ref{lem.sol.transpose.code.tij1} \textbf{(b)} (applied to $k$ and
$\sigma\left(  k\right)  $ instead of $i$ and $j$) yields $t_{k,\sigma\left(
k\right)  }\left(  \sigma\left(  k\right)  \right)  =k$. Thus, $\left(
t_{k,\sigma\left(  k\right)  }\circ\sigma\right)  \left(  k\right)
=t_{k,\sigma\left(  k\right)  }\left(  \sigma\left(  k\right)  \right)  =k=i$.
Hence, $\left(  t_{k,\sigma\left(  k\right)  }\circ\sigma\right)  \left(
\underbrace{i}_{=k}\right)  =\left(  t_{k,\sigma\left(  k\right)  }\circ
\sigma\right)  \left(  k\right)  =i$. Thus, $\left(  t_{k,\sigma\left(
k\right)  }\circ\sigma\right)  \left(  i\right)  =i$ is proven in Case 2.

We now have proven $\left(  t_{k,\sigma\left(  k\right)  }\circ\sigma\right)
\left(  i\right)  =i$ in each of the two Cases 1 and 2. Hence, $\left(
t_{k,\sigma\left(  k\right)  }\circ\sigma\right)  \left(  i\right)  =i$ always
holds. Thus, Lemma \ref{lem.sol.transpose.code.tij3} \textbf{(b)} is proven.

\textbf{(c)} Let $u\in\left[  n\right]  $ be such that $\left(  t_{k,u}%
\circ\sigma\right)  \left(  i\right)  =i$ for each $i\in\left\{
k,k+1,\ldots,n\right\}  $. We must prove that $u=\sigma\left(  k\right)  $.

We have assumed that $\left(  t_{k,u}\circ\sigma\right)  \left(  i\right)  =i$
for each $i\in\left\{  k,k+1,\ldots,n\right\}  $. Applying this to $i=k$, we
obtain $\left(  t_{k,u}\circ\sigma\right)  \left(  k\right)  =k$ (since
$k\in\left\{  k,k+1,\ldots,n\right\}  $). Hence, $t_{k,u}\left(  \sigma\left(
k\right)  \right)  =\left(  t_{k,u}\circ\sigma\right)  \left(  k\right)  =k$.

But Lemma \ref{lem.sol.transpose.code.tij1} \textbf{(b)} (applied to $k$ and
$u$ instead of $i$ and $j$) yields $t_{k,u}\left(  u\right)  =k$. Comparing
this with $t_{k,u}\left(  \sigma\left(  k\right)  \right)  =k$, we obtain
$t_{k,u}\left(  u\right)  =t_{k,u}\left(  \sigma\left(  k\right)  \right)  $.

But $t_{k,u}$ is an element of $S_{n}$. In other words, $t_{k,u}$ is a
permutation of $\left[  n\right]  $ (since $S_{n}$ is the set of all
permutations of $\left[  n\right]  $). Hence, the map $t_{k,u}$ is bijective,
and thus injective. From $t_{k,u}\left(  u\right)  =t_{k,u}\left(
\sigma\left(  k\right)  \right)  $, we therefore obtain $u=\sigma\left(
k\right)  $. This proves Lemma \ref{lem.sol.transpose.code.tij3} \textbf{(c)}.
\end{proof}

VERLONG PROOF:

\begin{proof}
[Proof of Lemma \ref{lem.sol.transpose.code.tij3}.]The definition of $\left[
k\right]  $ yields $\left[  k\right]  =\left\{  1,2,\ldots,k\right\}  $.

We have $k\in\left[  n\right]  =\left\{  1,2,\ldots,n\right\}  $ (by the
definition of $\left[  n\right]  $). Hence, $k$ is a positive integer. Thus,
$k\in\left\{  1,2,\ldots,k\right\}  =\left[  k\right]  $.

We have $\sigma\in S_{n}$. In other words, $\sigma$ is a permutation of
$\left\{  1,2,\ldots,n\right\}  $ (since $S_{n}$ is the set of all
permutations of $\left\{  1,2,\ldots,n\right\}  $). In other words, $\sigma$
is a permutation of $\left[  n\right]  $ (since $\left[  n\right]  =\left\{
1,2,\ldots,n\right\}  $). In other words, $\sigma$ is a bijection $\left[
n\right]  \rightarrow\left[  n\right]  $. Hence, $\sigma\left(  k\right)
\in\left[  n\right]  $.

\textbf{(a)} Assume the contrary. Thus, $\sigma\left(  k\right)  \notin\left[
k\right]  $.

Combining $\sigma\left(  k\right)  \in\left[  n\right]  =\left\{
1,2,\ldots,n\right\}  $ with $\sigma\left(  k\right)  \notin\left[  k\right]
=\left\{  1,2,\ldots,k\right\}  $, we obtain
\[
\sigma\left(  k\right)  \in\left\{  1,2,\ldots,n\right\}  \setminus\left\{
1,2,\ldots,k\right\}  =\left\{  k+1,k+2,\ldots,n\right\}  .
\]
Hence, (\ref{eq.lem.sol.transpose.code.tij3.ass}) (applied to $i=\sigma\left(
k\right)  $) yields $\sigma\left(  \sigma\left(  k\right)  \right)
=\sigma\left(  k\right)  $.

But the map $\sigma$ is bijective (since $\sigma$ is a bijection), thus
injective. Therefore, from $\sigma\left(  \sigma\left(  k\right)  \right)
=\sigma\left(  k\right)  $, we obtain $\sigma\left(  k\right)  =k\in\left[
k\right]  $. This contradicts $\sigma\left(  k\right)  \notin\left[  k\right]
$.

This contradiction shows that our assumption was wrong. Hence, Lemma
\ref{lem.sol.transpose.code.tij3} \textbf{(a)} is proven.

\textbf{(b)} Let $i\in\left\{  k,k+1,\ldots,n\right\}  $. We must prove that
$\left(  t_{k,\sigma\left(  k\right)  }\circ\sigma\right)  \left(  i\right)
=i$.

We are in one of the following two cases:

\textit{Case 1:} We have $i\neq k$.

\textit{Case 2:} We have $i=k$.

Let us first consider Case 1. In this case, we have $i\neq k$. Combining
$i\in\left\{  k,k+1,\ldots,n\right\}  $ with $i\neq k$, we obtain
$i\in\left\{  k,k+1,\ldots,n\right\}  \setminus\left\{  k\right\}  =\left\{
k+1,k+2,\ldots,n\right\}  $. Hence, (\ref{eq.lem.sol.transpose.code.tij3.ass})
shows that $\sigma\left(  i\right)  =i$.

But $i\notin\left\{  k,\sigma\left(  k\right)  \right\}  $%
\ \ \ \ \footnote{\textit{Proof.} Assume the contrary. Thus, $i\in\left\{
k,\sigma\left(  k\right)  \right\}  $. Combining this with $i\neq k$, we
obtain $i\in\left\{  k,\sigma\left(  k\right)  \right\}  \setminus\left\{
k\right\}  \subseteq\left\{  \sigma\left(  k\right)  \right\}  $. Thus,
$i=\sigma\left(  k\right)  \in\left[  k\right]  $ (by Lemma
\ref{lem.sol.transpose.code.tij3} \textbf{(a)}). Hence, $i\in\left[  k\right]
=\left\{  1,2,\ldots,k\right\}  $, so that $i\leq k$. But $i\in\left\{
k+1,k+2,\ldots,n\right\}  $, so that $i\geq k+1>k$. This contradicts $i\leq
k$. This contradiction shows that our assumption was wrong. Qed.}. Combining
$i\in\left\{  k,k+1,\ldots,n\right\}  \subseteq\left\{  1,2,\ldots,n\right\}
=\left[  n\right]  $ with $i\notin\left\{  k,\sigma\left(  k\right)  \right\}
$, we obtain $i\in\left[  n\right]  \setminus\left\{  k,\sigma\left(
k\right)  \right\}  $. Hence, Lemma \ref{lem.sol.transpose.code.tij1}
\textbf{(c)} (applied to $k$, $\sigma\left(  k\right)  $ and $i$ instead of
$i$, $j$ and $k$) shows that $t_{k,\sigma\left(  k\right)  }\left(  i\right)
=i$.

Now, $\left(  t_{k,\sigma\left(  k\right)  }\circ\sigma\right)  \left(
i\right)  =t_{k,\sigma\left(  k\right)  }\left(  \underbrace{\sigma\left(
i\right)  }_{=i}\right)  =t_{k,\sigma\left(  k\right)  }\left(  i\right)  =i$.
Thus, $\left(  t_{k,\sigma\left(  k\right)  }\circ\sigma\right)  \left(
i\right)  =i$ is proven in Case 1.

Let us now consider Case 2. In this case, we have $i=k$. On the other hand,
Lemma \ref{lem.sol.transpose.code.tij1} \textbf{(b)} (applied to $k$ and
$\sigma\left(  k\right)  $ instead of $i$ and $j$) yields $t_{k,\sigma\left(
k\right)  }\left(  \sigma\left(  k\right)  \right)  =k$. Thus, $\left(
t_{k,\sigma\left(  k\right)  }\circ\sigma\right)  \left(  k\right)
=t_{k,\sigma\left(  k\right)  }\left(  \sigma\left(  k\right)  \right)  =k=i$.
Hence, $\left(  t_{k,\sigma\left(  k\right)  }\circ\sigma\right)  \left(
\underbrace{i}_{=k}\right)  =\left(  t_{k,\sigma\left(  k\right)  }\circ
\sigma\right)  \left(  k\right)  =i$. Thus, $\left(  t_{k,\sigma\left(
k\right)  }\circ\sigma\right)  \left(  i\right)  =i$ is proven in Case 2.

We now have proven $\left(  t_{k,\sigma\left(  k\right)  }\circ\sigma\right)
\left(  i\right)  =i$ in each of the two Cases 1 and 2. Since these two Cases
cover all possibilities, we thus conclude that $\left(  t_{k,\sigma\left(
k\right)  }\circ\sigma\right)  \left(  i\right)  =i$ always holds. Thus, Lemma
\ref{lem.sol.transpose.code.tij3} \textbf{(b)} is proven.

\textbf{(c)} Let $u\in\left[  n\right]  $ be such that $\left(  t_{k,u}%
\circ\sigma\right)  \left(  i\right)  =i$ for each $i\in\left\{
k,k+1,\ldots,n\right\}  $. We must prove that $u=\sigma\left(  k\right)  $.

We have $k\in\left\{  1,2,\ldots,n\right\}  $, so that $k\leq n$ and thus
$k\leq k\leq n$. Hence, $k\in\left\{  k,k+1,\ldots,n\right\}  $. But we have
assumed that $\left(  t_{k,u}\circ\sigma\right)  \left(  i\right)  =i$ for
each $i\in\left\{  k,k+1,\ldots,n\right\}  $. Applying this to $i=k$, we
obtain $\left(  t_{k,u}\circ\sigma\right)  \left(  k\right)  =k$ (since
$k\in\left\{  k,k+1,\ldots,n\right\}  $). Hence, $t_{k,u}\left(  \sigma\left(
k\right)  \right)  =\left(  t_{k,u}\circ\sigma\right)  \left(  k\right)  =k$.

But Lemma \ref{lem.sol.transpose.code.tij1} \textbf{(b)} (applied to $k$ and
$u$ instead of $i$ and $j$) yields $t_{k,u}\left(  u\right)  =k$. Comparing
this with $t_{k,u}\left(  \sigma\left(  k\right)  \right)  =k$, we obtain
$t_{k,u}\left(  \sigma\left(  k\right)  \right)  =t_{k,u}\left(  u\right)  $.

But $t_{k,u}$ is an element of $S_{n}$. In other words, $t_{k,u}$ is a
permutation of $\left\{  1,2,\ldots,n\right\}  $ (since $S_{n}$ is the set of
all permutations of $\left\{  1,2,\ldots,n\right\}  $). In other words,
$t_{k,u}$ is a bijection $\left\{  1,2,\ldots,n\right\}  \rightarrow\left\{
1,2,\ldots,n\right\}  $. Hence, the map $t_{k,u}$ is bijective, and thus
injective. From $t_{k,u}\left(  \sigma\left(  k\right)  \right)
=t_{k,u}\left(  u\right)  $, we therefore obtain $\sigma\left(  k\right)  =u$.
In other words, $u=\sigma\left(  k\right)  $. This proves Lemma
\ref{lem.sol.transpose.code.tij3} \textbf{(c)}.
\end{proof}
\end{noncompile}

\begin{lemma}
\label{lem.sol.transpose.code.exist}Let $n\in\mathbb{N}$. Let $k\in\left\{
0,1,\ldots,n\right\}  $ and $\sigma\in S_{n}$. Assume that%
\begin{equation}
\left(  \sigma\left(  i\right)  =i\text{ for each }i\in\left\{  k+1,k+2,\ldots
,n\right\}  \right)  . \label{eq.lem.sol.transpose.code.exist.ass}%
\end{equation}
Then, there exists a $k$-tuple $\left(  i_{1},i_{2},\ldots,i_{k}\right)
\in\left[  1\right]  \times\left[  2\right]  \times\cdots\times\left[
k\right]  $ such that%
\[
\sigma=t_{1,i_{1}}\circ t_{2,i_{2}}\circ\cdots\circ t_{k,i_{k}}.
\]

\end{lemma}

\begin{proof}
[Proof of Lemma \ref{lem.sol.transpose.code.exist}.]We shall prove Lemma
\ref{lem.sol.transpose.code.exist} by induction over $k$:

\begin{vershort}
\textit{Induction base:} Lemma \ref{lem.sol.transpose.code.exist} is true when
$k=0$\ \ \ \ \footnote{\textit{Proof.} Let $n$, $k$ and $\sigma$ be as in
Lemma \ref{lem.sol.transpose.code.exist}. Assume that $k=0$. We have to prove
that Lemma \ref{lem.sol.transpose.code.exist} holds under this assumption.
\par
In view of $k=0$, the assumption (\ref{eq.lem.sol.transpose.code.exist.ass})
rewrites as follows: We have $\sigma\left(  i\right)  =i$ for each
$i\in\left\{  1,2,\ldots,n\right\}  $. In other words, $\sigma
=\operatorname*{id}$.
\par
If $\left(  i_{1},i_{2},\ldots,i_{0}\right)  \in\left[  1\right]
\times\left[  2\right]  \times\cdots\times\left[  0\right]  $ is the empty
$0$-tuple $\left(  {}\right)  $, then%
\[
t_{1,i_{1}}\circ t_{2,i_{2}}\circ\cdots\circ t_{0,i_{0}}=\left(  \text{empty
composition of permutations}\right)  =\operatorname*{id}=\sigma.
\]
Hence, there exists a $0$-tuple $\left(  i_{1},i_{2},\ldots,i_{0}\right)
\in\left[  1\right]  \times\left[  2\right]  \times\cdots\times\left[
0\right]  $ such that $\sigma=t_{1,i_{1}}\circ t_{2,i_{2}}\circ\cdots\circ
t_{0,i_{0}}$ (namely, the empty $0$-tuple $\left(  {}\right)  $). Since $k=0$,
this rewrites as follows: There exists a $k$-tuple $\left(  i_{1},i_{2}%
,\ldots,i_{k}\right)  \in\left[  1\right]  \times\left[  2\right]
\times\cdots\times\left[  k\right]  $ such that $\sigma=t_{1,i_{1}}\circ
t_{2,i_{2}}\circ\cdots\circ t_{k,i_{k}}$. In other words, Lemma
\ref{lem.sol.transpose.code.exist} holds. Thus, Lemma
\ref{lem.sol.transpose.code.exist} is true when $k=0$.}.
\end{vershort}

\begin{verlong}
\textit{Induction base:} Lemma \ref{lem.sol.transpose.code.exist} is true when
$k=0$\ \ \ \ \footnote{\textit{Proof.} Let $n$, $k$ and $\sigma$ be as in
Lemma \ref{lem.sol.transpose.code.exist}. Assume that $k=0$. We have to prove
that Lemma \ref{lem.sol.transpose.code.exist} holds under this assumption.
\par
Let $i\in\left\{  1,2,\ldots,n\right\}  $. Notice that $k=0$, so that $k+1=1$
and therefore $\left\{  k+1,k+2,\ldots,n\right\}  =\left\{  1,2,\ldots
,n\right\}  $. Now, $i\in\left\{  1,2,\ldots,n\right\}  =\left\{
k+1,k+2,\ldots,n\right\}  $. Hence, (\ref{eq.lem.sol.transpose.code.exist.ass}%
) shows that $\sigma\left(  i\right)  =i=\operatorname*{id}\left(  i\right)
$.
\par
Now, forget that we fixed $i$. We thus have proven that $\sigma\left(
i\right)  =\operatorname*{id}\left(  i\right)  $ for each $i\in\left\{
1,2,\ldots,n\right\}  $. In other words, $\sigma=\operatorname*{id}$ (since
both $\sigma$ and $\operatorname*{id}$ are maps $\left\{  1,2,\ldots
,n\right\}  \rightarrow\left\{  1,2,\ldots,n\right\}  $).
\par
If $\left(  i_{1},i_{2},\ldots,i_{0}\right)  \in\left[  1\right]
\times\left[  2\right]  \times\cdots\times\left[  0\right]  $ is the empty
$0$-tuple $\left(  {}\right)  $, then
\[
t_{1,i_{1}}\circ t_{2,i_{2}}\circ\cdots\circ t_{0,i_{0}}=\left(  \text{empty
composition of permutations}\right)  =\operatorname*{id}=\sigma
\]
(since $\sigma=\operatorname*{id}$). In other words, if $\left(  i_{1}%
,i_{2},\ldots,i_{0}\right)  \in\left[  1\right]  \times\left[  2\right]
\times\cdots\times\left[  0\right]  $ is the empty $0$-tuple $\left(
{}\right)  $, then $\sigma=t_{1,i_{1}}\circ t_{2,i_{2}}\circ\cdots\circ
t_{0,i_{0}}$. Hence, there exists a $0$-tuple $\left(  i_{1},i_{2}%
,\ldots,i_{0}\right)  \in\left[  1\right]  \times\left[  2\right]
\times\cdots\times\left[  0\right]  $ such that $\sigma=t_{1,i_{1}}\circ
t_{2,i_{2}}\circ\cdots\circ t_{0,i_{0}}$ (namely, the empty $0$-tuple $\left(
{}\right)  $). Since $k=0$, this rewrites as follows: There exists a $k$-tuple
$\left(  i_{1},i_{2},\ldots,i_{k}\right)  \in\left[  1\right]  \times\left[
2\right]  \times\cdots\times\left[  k\right]  $ such that $\sigma=t_{1,i_{1}%
}\circ t_{2,i_{2}}\circ\cdots\circ t_{k,i_{k}}$.
\par
Now, forget that we fixed $n$, $k$ and $\sigma$. We thus have proven that if
$n$, $k$ and $\sigma$ are as in Lemma \ref{lem.sol.transpose.code.exist}, and
if we have $k=0$, then there exists a $k$-tuple $\left(  i_{1},i_{2}%
,\ldots,i_{k}\right)  \in\left[  1\right]  \times\left[  2\right]
\times\cdots\times\left[  k\right]  $ such that $\sigma=t_{1,i_{1}}\circ
t_{2,i_{2}}\circ\cdots\circ t_{k,i_{k}}$. In other words, Lemma
\ref{lem.sol.transpose.code.exist} is true when $k=0$. Qed.}. This completes
the induction base.
\end{verlong}

\textit{Induction step:} Let $K\in\left\{  0,1,\ldots,n\right\}  $ be
positive. Assume that Lemma \ref{lem.sol.transpose.code.exist} holds when
$k=K-1$. We must show that Lemma \ref{lem.sol.transpose.code.exist} holds when
$k=K$.

We have assumed that Lemma \ref{lem.sol.transpose.code.exist} holds when
$k=K-1$. In other words, the following statement holds:

\begin{statement}
\textit{Statement 1:} Let $n\in\mathbb{N}$. Assume that $K-1\in\left\{
0,1,\ldots,n\right\}  $. Let $\sigma\in S_{n}$. Assume that%
\[
\left(  \sigma\left(  i\right)  =i\text{ for each }i\in\left\{  \left(
K-1\right)  +1,\left(  K-1\right)  +2,\ldots,n\right\}  \right)  .
\]
Then, there exists a $\left(  K-1\right)  $-tuple $\left(  i_{1},i_{2}%
,\ldots,i_{K-1}\right)  \in\left[  1\right]  \times\left[  2\right]
\times\cdots\times\left[  K-1\right]  $ such that%
\[
\sigma=t_{1,i_{1}}\circ t_{2,i_{2}}\circ\cdots\circ t_{K-1,i_{K-1}}.
\]

\end{statement}

We must show that Lemma \ref{lem.sol.transpose.code.exist} holds when $k=K$.
In other words, we must prove the following statement:

\begin{statement}
\textit{Statement 2:} Let $n\in\mathbb{N}$. Assume that $K\in\left\{
0,1,\ldots,n\right\}  $. Let $\sigma\in S_{n}$. Assume that%
\[
\left(  \sigma\left(  i\right)  =i\text{ for each }i\in\left\{  K+1,K+2,\ldots
,n\right\}  \right)  .
\]
Then, there exists a $K$-tuple $\left(  i_{1},i_{2},\ldots,i_{K}\right)
\in\left[  1\right]  \times\left[  2\right]  \times\cdots\times\left[
K\right]  $ such that%
\[
\sigma=t_{1,i_{1}}\circ t_{2,i_{2}}\circ\cdots\circ t_{K,i_{K}}.
\]

\end{statement}

\begin{vershort}
[\textit{Proof of Statement 2:} We have $K\neq0$ (since $K$ is positive).
Combined with $K\in\left\{  0,1,\ldots,n\right\}  $, this yields $K\in\left\{
0,1,\ldots,n\right\}  \setminus\left\{  0\right\}  =\left\{  1,2,\ldots
,n\right\}  $, so that $K-1\in\left\{  0,1,\ldots,n-1\right\}  \subseteq
\left\{  0,1,\ldots,n\right\}  $.

We have $\sigma\in S_{n}$. In other words, $\sigma$ is a permutation of
$\left[  n\right]  $ (since $S_{n}$ is the set of all permutations of $\left[
n\right]  $). Hence, $\sigma$ has an inverse $\sigma^{-1}$.

Moreover, $K\in\left\{  1,2,\ldots,n\right\}  =\left[  n\right]  $. We can
thus define $g\in\left[  n\right]  $ by $g=\sigma^{-1}\left(  K\right)  $.
Consider this $g$.

Lemma \ref{lem.sol.transpose.code.tij2} \textbf{(a)} (applied to $k=K$) yields
$g\in\left[  K\right]  $.

Define a permutation $\tau\in S_{n}$ by $\tau=\sigma\circ t_{K,g}$. Lemma
\ref{lem.sol.transpose.code.tij2} \textbf{(b)} (applied to $k=K$) yields that
$\left(  \sigma\circ t_{K,g}\right)  \left(  i\right)  =i$ for each
$i\in\left\{  K,K+1,\ldots,n\right\}  $. In other words, $\tau\left(
i\right)  =i$ for each $i\in\left\{  K,K+1,\ldots,n\right\}  $ (since
$\tau=\sigma\circ t_{K,g}$). In other words, $\tau\left(  i\right)  =i$ for
each $i\in\left\{  \left(  K-1\right)  +1,\left(  K-1\right)  +2,\ldots
,n\right\}  $ (since $K=\left(  K-1\right)  +1$). Hence, Statement 1 (applied
to $\tau$ instead of $\sigma$) shows that there exists a $\left(  K-1\right)
$-tuple $\left(  i_{1},i_{2},\ldots,i_{K-1}\right)  \in\left[  1\right]
\times\left[  2\right]  \times\cdots\times\left[  K-1\right]  $ such that%
\[
\tau=t_{1,i_{1}}\circ t_{2,i_{2}}\circ\cdots\circ t_{K-1,i_{K-1}}.
\]
Consider this $\left(  K-1\right)  $-tuple $\left(  i_{1},i_{2},\ldots
,i_{K-1}\right)  $, and denote it by $\left(  j_{1},j_{2},\ldots
,j_{K-1}\right)  $. Thus, $\left(  j_{1},j_{2},\ldots,j_{K-1}\right)  $ is a
$\left(  K-1\right)  $-tuple in $\left[  1\right]  \times\left[  2\right]
\times\cdots\times\left[  K-1\right]  $ such that%
\[
\tau=t_{1,j_{1}}\circ t_{2,j_{2}}\circ\cdots\circ t_{K-1,j_{K-1}}.
\]
Now,
\[
\underbrace{\tau}_{=\sigma\circ t_{K,g}}\circ t_{K,g}=\sigma\circ
\underbrace{t_{K,g}\circ t_{K,g}}_{\substack{=\operatorname*{id}\\\text{(by
Lemma \ref{lem.sol.transpose.code.tij1} \textbf{(d)},}\\\text{applied to
}i=K\text{ and }j=g\text{)}}}=\sigma.
\]
Hence,%
\begin{equation}
\sigma=\underbrace{\tau}_{=t_{1,j_{1}}\circ t_{2,j_{2}}\circ\cdots\circ
t_{K-1,j_{K-1}}}\circ t_{K,g}=\left(  t_{1,j_{1}}\circ t_{2,j_{2}}\circ
\cdots\circ t_{K-1,j_{K-1}}\right)  \circ t_{K,g}.
\label{pf.lem.sol.transpose.code.exist.st2.pf.3}%
\end{equation}

Now, let us extend our $\left(  K-1\right)  $-tuple $\left(  j_{1}%
,j_{2},\ldots,j_{K-1}\right)  \in\left[  1\right]  \times\left[  2\right]
\times\cdots\times\left[  K-1\right]  $ to a $K$-tuple $\left(  j_{1}%
,j_{2},\ldots,j_{K}\right)  \in\left[  1\right]  \times\left[  2\right]
\times\cdots\times\left[  K\right]  $ by setting $j_{K}=g$. (This is
well-defined, since $g\in\left[  K\right]  $.) Then,
\begin{align*}
t_{1,j_{1}}\circ t_{2,j_{2}}\circ\cdots\circ t_{K,j_{K}}  &  =\left(
t_{1,j_{1}}\circ t_{2,j_{2}}\circ\cdots\circ t_{K-1,j_{K-1}}\right)
\circ\underbrace{t_{K,j_{K}}}_{\substack{=t_{K,g}\\\text{(since }%
j_{K}=g\text{)}}}\\
&  =\left(  t_{1,j_{1}}\circ t_{2,j_{2}}\circ\cdots\circ t_{K-1,j_{K-1}%
}\right)  \circ t_{K,g}.
\end{align*}
Comparing this with (\ref{pf.lem.sol.transpose.code.exist.st2.pf.3}), we find
$\sigma=t_{1,j_{1}}\circ t_{2,j_{2}}\circ\cdots\circ t_{K,j_{K}}$. Thus, there
exists a $K$-tuple $\left(  i_{1},i_{2},\ldots,i_{K}\right)  \in\left[
1\right]  \times\left[  2\right]  \times\cdots\times\left[  K\right]  $ such
that%
\[
\sigma=t_{1,i_{1}}\circ t_{2,i_{2}}\circ\cdots\circ t_{K,i_{K}}%
\]
(namely, $\left(  i_{1},i_{2},\ldots,i_{K}\right)  =\left(  j_{1},j_{2}%
,\ldots,j_{K}\right)  $). This proves Statement 2.]
\end{vershort}

\begin{verlong}
[\textit{Proof of Statement 2:} We have $\left[  n\right]  =\left\{
1,2,\ldots,n\right\}  $ (by the definition of $\left[  n\right]  $).

We have $\left\{  K,K+1,\ldots,n\right\}  =\left\{  \left(  K-1\right)
+1,\left(  K-1\right)  +2,\ldots,n\right\}  $ (since $K=\left(  K-1\right)
+1$).

We have $K\neq0$ (since $K$ is positive). Combining $K\in\left\{
0,1,\ldots,n\right\}  $ with $K\neq0$, we obtain $K\in\left\{  0,1,\ldots
,n\right\}  \setminus\left\{  0\right\}  =\left\{  1,2,\ldots,n\right\}  $, so
that $K-1\in\left\{  0,1,\ldots,n-1\right\}  \subseteq\left\{  0,1,\ldots
,n\right\}  $.

We have $\sigma\in S_{n}$. In other words, $\sigma$ is a permutation of
$\left\{  1,2,\ldots,n\right\}  $ (since $S_{n}$ is the set of all
permutations of $\left\{  1,2,\ldots,n\right\}  $). In other words, $\sigma$
is a permutation of $\left[  n\right]  $ (since $\left[  n\right]  =\left\{
1,2,\ldots,n\right\}  $). In other words, $\sigma$ is a bijection $\left[
n\right]  \rightarrow\left[  n\right]  $. Hence, this map $\sigma$ is
invertible, and its inverse $\sigma^{-1}$ is also a bijection $\left[
n\right]  \rightarrow\left[  n\right]  $.

Moreover, $K\in\left\{  1,2,\ldots,n\right\}  =\left[  n\right]  $. Hence,
$\sigma^{-1}\left(  K\right)  \in\left[  n\right]  $ (since $\sigma^{-1}$ is a
map $\left[  n\right]  \rightarrow\left[  n\right]  $). We can thus define
$g\in\left[  n\right]  $ by $g=\sigma^{-1}\left(  K\right)  $. Consider this
$g$.

Lemma \ref{lem.sol.transpose.code.tij2} \textbf{(a)} (applied to $k=K$) yields
$g\in\left[  K\right]  $.

Define a permutation $\tau\in S_{n}$ by $\tau=\sigma\circ t_{K,g}$. Lemma
\ref{lem.sol.transpose.code.tij2} \textbf{(b)} (applied to $k=K$) yields that
$\left(  \sigma\circ t_{K,g}\right)  \left(  i\right)  =i$ for each
$i\in\left\{  K,K+1,\ldots,n\right\}  $. In other words, $\tau\left(
i\right)  =i$ for each $i\in\left\{  K,K+1,\ldots,n\right\}  $ (since
$\tau=\sigma\circ t_{K,g}$). In other words, $\tau\left(  i\right)  =i$ for
each $i\in\left\{  \left(  K-1\right)  +1,\left(  K-1\right)  +2,\ldots
,n\right\}  $ (since $\left\{  K,K+1,\ldots,n\right\}  =\left\{  \left(
K-1\right)  +1,\left(  K-1\right)  +2,\ldots,n\right\}  $). Hence, Statement 1
(applied to $\tau$ instead of $\sigma$) shows that there exists a $\left(
K-1\right)  $-tuple $\left(  i_{1},i_{2},\ldots,i_{K-1}\right)  \in\left[
1\right]  \times\left[  2\right]  \times\cdots\times\left[  K-1\right]  $ such
that%
\[
\tau=t_{1,i_{1}}\circ t_{2,i_{2}}\circ\cdots\circ t_{K-1,i_{K-1}}.
\]
Consider this $\left(  K-1\right)  $-tuple $\left(  i_{1},i_{2},\ldots
,i_{K-1}\right)  $, and denote it by $\left(  j_{1},j_{2},\ldots
,j_{K-1}\right)  $. Thus, $\left(  j_{1},j_{2},\ldots,j_{K-1}\right)  $ is a
$\left(  K-1\right)  $-tuple in $\left[  1\right]  \times\left[  2\right]
\times\cdots\times\left[  K-1\right]  $ such that%
\[
\tau=t_{1,j_{1}}\circ t_{2,j_{2}}\circ\cdots\circ t_{K-1,j_{K-1}}.
\]
Now, Lemma \ref{lem.sol.transpose.code.tij1} \textbf{(d)} (applied to $K$ and
$g$ instead of $i$ and $j$) shows that $t_{K,g}\circ t_{K,g}%
=\operatorname*{id}$. Now,
\[
\underbrace{\tau}_{=\sigma\circ t_{K,g}}\circ t_{K,g}=\sigma\circ
\underbrace{t_{K,g}\circ t_{K,g}}_{=\operatorname*{id}}=\sigma.
\]
Hence,%
\[
\sigma=\underbrace{\tau}_{=t_{1,j_{1}}\circ t_{2,j_{2}}\circ\cdots\circ
t_{K-1,j_{K-1}}}\circ t_{K,g}=\left(  t_{1,j_{1}}\circ t_{2,j_{2}}\circ
\cdots\circ t_{K-1,j_{K-1}}\right)  \circ t_{K,g}.
\]
In other words, $\left(  t_{1,j_{1}}\circ t_{2,j_{2}}\circ\cdots\circ
t_{K-1,j_{K-1}}\right)  \circ t_{K,g}=\sigma$.

Now, let us extend our $\left(  K-1\right)  $-tuple $\left(  j_{1}%
,j_{2},\ldots,j_{K-1}\right)  \in\left[  1\right]  \times\left[  2\right]
\times\cdots\times\left[  K-1\right]  $ to a $K$-tuple $\left(  j_{1}%
,j_{2},\ldots,j_{K}\right)  \in\left[  1\right]  \times\left[  2\right]
\times\cdots\times\left[  K\right]  $ by setting $j_{K}=g$. (This is
well-defined, since $g\in\left[  K\right]  $.) Then,
\begin{align*}
t_{1,j_{1}}\circ t_{2,j_{2}}\circ\cdots\circ t_{K,j_{K}}  &  =\left(
t_{1,j_{1}}\circ t_{2,j_{2}}\circ\cdots\circ t_{K-1,j_{K-1}}\right)
\circ\underbrace{t_{K,j_{K}}}_{\substack{=t_{K,g}\\\text{(since }%
j_{K}=g\text{)}}}\\
&  =\left(  t_{1,j_{1}}\circ t_{2,j_{2}}\circ\cdots\circ t_{K-1,j_{K-1}%
}\right)  \circ t_{K,g}=\sigma.
\end{align*}
Hence, $\sigma=t_{1,j_{1}}\circ t_{2,j_{2}}\circ\cdots\circ t_{K,j_{K}}$.
Thus, there exists a $K$-tuple $\left(  i_{1},i_{2},\ldots,i_{K}\right)
\in\left[  1\right]  \times\left[  2\right]  \times\cdots\times\left[
K\right]  $ such that%
\[
\sigma=t_{1,i_{1}}\circ t_{2,i_{2}}\circ\cdots\circ t_{K,i_{K}}%
\]
(namely, $\left(  i_{1},i_{2},\ldots,i_{K}\right)  =\left(  j_{1},j_{2}%
,\ldots,j_{K}\right)  $). This proves Statement 2.]
\end{verlong}

So Statement 2 is proven. In other words, Lemma
\ref{lem.sol.transpose.code.exist} holds when $k=K$ (since Statement 2 is
precisely the claim of Lemma \ref{lem.sol.transpose.code.exist} for $k=K$).
This completes the induction step. Thus, Lemma
\ref{lem.sol.transpose.code.exist} is proven.
\end{proof}

\begin{lemma}
\label{lem.sol.transpose.code.uni}Let $n\in\mathbb{N}$. Let $k\in\left\{
0,1,\ldots,n\right\}  $ and $\sigma\in S_{n}$. Let $\left(  i_{1},i_{2}%
,\ldots,i_{k}\right)  \in\left[  1\right]  \times\left[  2\right]
\times\cdots\times\left[  k\right]  $ be such that $\sigma=t_{1,i_{1}}\circ
t_{2,i_{2}}\circ\cdots\circ t_{k,i_{k}}$. Then:

\textbf{(a)} We have $\sigma\left(  i\right)  =i$ for each $i\in\left\{
k+1,k+2,\ldots,n\right\}  $.

\textbf{(b)} If $k>0$, then $\sigma\left(  i_{k}\right)  =k$.
\end{lemma}

\begin{vershort}
\begin{proof}
[Proof of Lemma \ref{lem.sol.transpose.code.uni}.]We have $\left(  i_{1}%
,i_{2},\ldots,i_{k}\right)  \in\left[  1\right]  \times\left[  2\right]
\times\cdots\times\left[  k\right]  $. Thus,%
\begin{equation}
i_{j}\in\left[  j\right]  \ \ \ \ \ \ \ \ \ \ \text{for each }j\in\left\{
1,2,\ldots,k\right\}  . \label{pf.lem.sol.transpose.code.uni.short.1}%
\end{equation}

\textbf{(a)} Let $i\in\left\{  k+1,k+2,\ldots,n\right\}  $. Thus, $i\geq
k+1>k$.

We claim that%
\begin{equation}
\left(  t_{1,i_{1}}\circ t_{2,i_{2}}\circ\cdots\circ t_{p,i_{p}}\right)
\left(  i\right)  =i\ \ \ \ \ \ \ \ \ \ \text{for each }p\in\left\{
0,1,\ldots,k\right\}  . \label{pf.lem.sol.transpose.code.uni.short.c1.1}%
\end{equation}

[\textit{Proof of (\ref{pf.lem.sol.transpose.code.uni.short.c1.1}):} We shall
prove (\ref{pf.lem.sol.transpose.code.uni.short.c1.1}) by induction over $p$:

\textit{Induction base:} We have $\underbrace{\left(  t_{1,i_{1}}\circ
t_{2,i_{2}}\circ\cdots\circ t_{0,i_{0}}\right)  }_{=\left(  \text{empty
composition}\right)  =\operatorname*{id}}\left(  i\right)  =\operatorname*{id}%
\left(  i\right)  =i$. In other words,
(\ref{pf.lem.sol.transpose.code.uni.short.c1.1}) holds for $p=0$. This
completes the induction base.

\textit{Induction step:} Let $q\in\left\{  0,1,\ldots,k\right\}  $ be
positive. Assume that (\ref{pf.lem.sol.transpose.code.uni.short.c1.1}) holds
for $p=q-1$. We must prove that
(\ref{pf.lem.sol.transpose.code.uni.short.c1.1}) holds for $p=q$.

We have assumed that (\ref{pf.lem.sol.transpose.code.uni.short.c1.1}) holds
for $p=q-1$. In other words, we have $\left(  t_{1,i_{1}}\circ t_{2,i_{2}%
}\circ\cdots\circ t_{q-1,i_{q-1}}\right)  \left(  i\right)  =i$.

We have $q\in\left\{  0,1,\ldots,k\right\}  $, so that $q\leq k$. Hence,
$k\geq q$, so that $i>k\geq q$. Thus, $i\neq q$. Furthermore, $q\neq0$ (since
$q$ is positive). Combining this with $q\in\left\{  0,1,\ldots,k\right\}  $,
we obtain $q\in\left\{  0,1,\ldots,k\right\}  \setminus\left\{  0\right\}
=\left\{  1,2,\ldots,k\right\}  $. Hence,
(\ref{pf.lem.sol.transpose.code.uni.short.1}) (applied to $j=q$) yields
$i_{q}\in\left[  q\right]  $. Hence, $i_{q}\leq q$. Thus, $q\geq i_{q}$, so
that $i>q\geq i_{q}$. Thus, $i\neq i_{q}$. Finally, $i_{q}\in\left[  q\right]
\subseteq\left[  n\right]  $ (since $q\leq k\leq n$). Furthermore, $q\leq
k\leq n$, so that $q\in\left[  n\right]  $ (since $q$ is positive).

Also, $i\in\left\{  k+1,k+2,\ldots,n\right\}  \subseteq\left\{  1,2,\ldots
,n\right\}  =\left[  n\right]  $. Thus, $i$ is an element of $\left[
n\right]  $ other than $q$ and $i_{q}$ (since $i\neq q$ and $i\neq i_{q}$). In
other words, $i\in\left[  n\right]  \setminus\left\{  q,i_{q}\right\}  $.
Hence, Lemma \ref{lem.sol.transpose.code.tij1} \textbf{(c)} (applied to $q$,
$i_{q}$ and $i$ instead of $i$, $j$ and $k$) shows that $t_{q,i_{q}}\left(
i\right)  =i$. Now,%
\begin{align*}
&  \underbrace{\left(  t_{1,i_{1}}\circ t_{2,i_{2}}\circ\cdots\circ
t_{q,i_{q}}\right)  }_{=\left(  t_{1,i_{1}}\circ t_{2,i_{2}}\circ\cdots\circ
t_{q-1,i_{q-1}}\right)  \circ t_{q,i_{q}}}\left(  i\right) \\
&  =\left(  \left(  t_{1,i_{1}}\circ t_{2,i_{2}}\circ\cdots\circ
t_{q-1,i_{q-1}}\right)  \circ t_{q,i_{q}}\right)  \left(  i\right)  =\left(
t_{1,i_{1}}\circ t_{2,i_{2}}\circ\cdots\circ t_{q-1,i_{q-1}}\right)  \left(
\underbrace{t_{q,i_{q}}\left(  i\right)  }_{=i}\right) \\
&  =\left(  t_{1,i_{1}}\circ t_{2,i_{2}}\circ\cdots\circ t_{q-1,i_{q-1}%
}\right)  \left(  i\right)  =i.
\end{align*}
In other words, (\ref{pf.lem.sol.transpose.code.uni.short.c1.1}) holds for
$p=q$. This completes the induction step. Thus,
(\ref{pf.lem.sol.transpose.code.uni.short.c1.1}) is proven.]

Now, we can apply (\ref{pf.lem.sol.transpose.code.uni.short.c1.1}) to $p=k$
(since $k\in\left\{  0,1,\ldots,k\right\}  $). We thus obtain $\left(
t_{1,i_{1}}\circ t_{2,i_{2}}\circ\cdots\circ t_{k,i_{k}}\right)  \left(
i\right)  =i$. Hence,%
\[
\underbrace{\sigma}_{=t_{1,i_{1}}\circ t_{2,i_{2}}\circ\cdots\circ t_{k,i_{k}%
}}\left(  i\right)  =\left(  t_{1,i_{1}}\circ t_{2,i_{2}}\circ\cdots\circ
t_{k,i_{k}}\right)  \left(  i\right)  =i.
\]
This proves Lemma \ref{lem.sol.transpose.code.uni} \textbf{(a)}.

\textbf{(b)} Assume that $k>0$. Hence, $k-1\in\mathbb{N}$. Hence,
$k-1\in\left\{  0,1,\ldots,n\right\}  $ (since $k-1\leq k\leq n$).
Furthermore, $k\in\left\{  \left(  k-1\right)  +1,\left(  k-1\right)
+2,\ldots,n\right\}  $ (because $k=\left(  k-1\right)  +1\geq\left(
k-1\right)  +1$ and $k\leq n$).

Define $\tau\in S_{n}$ by $\tau=t_{1,i_{1}}\circ t_{2,i_{2}}\circ\cdots\circ
t_{k-1,i_{k-1}}$. From (\ref{pf.lem.sol.transpose.code.uni.short.1}), we
obtain $\left(  i_{1},i_{2},\ldots,i_{k-1}\right)  \in\left[  1\right]
\times\left[  2\right]  \times\cdots\times\left[  k-1\right]  $. Hence, Lemma
\ref{lem.sol.transpose.code.uni} \textbf{(a)} (applied to $k-1$, $\tau$ and
$k$ instead of $k$, $\sigma$ and $i$) yields $\tau\left(  k\right)  =k$.

Combining $k>0$ with $k\in\left\{  0,1,\ldots,n\right\}  $, we obtain
$k\in\left[  n\right]  $. Also, $k\in\left\{  1,2,\ldots,k\right\}  $ (since
$k>0$). Hence, (\ref{pf.lem.sol.transpose.code.uni.short.1}) (applied to
$j=k$) yields $i_{k}\in\left[  k\right]  \subseteq\left[  n\right]  $ (since
$k\leq n$). Thus, Lemma \ref{lem.sol.transpose.code.tij1} \textbf{(b)}
(applied to $k$ and $i_{k}$ instead of $i$ and $j$) shows that $t_{k,i_{k}%
}\left(  i_{k}\right)  =k$. Now,%
\[
\sigma=t_{1,i_{1}}\circ t_{2,i_{2}}\circ\cdots\circ t_{k,i_{k}}%
=\underbrace{\left(  t_{1,i_{1}}\circ t_{2,i_{2}}\circ\cdots\circ
t_{k-1,i_{k-1}}\right)  }_{=\tau}\circ t_{k,i_{k}}=\tau\circ t_{k,i_{k}}.
\]
Applying both sides of this equality to $i_{k}$, we find%
\[
\sigma\left(  i_{k}\right)  =\left(  \tau\circ t_{k,i_{k}}\right)  \left(
i_{k}\right)  =\tau\left(  \underbrace{t_{k,i_{k}}\left(  i_{k}\right)  }%
_{=k}\right)  =\tau\left(  k\right)  =k.
\]
This proves Lemma \ref{lem.sol.transpose.code.uni} \textbf{(b)}.
\end{proof}
\end{vershort}

\begin{verlong}
\begin{proof}
[Proof of Lemma \ref{lem.sol.transpose.code.uni}.]We have $k\in\left\{
0,1,\ldots,n\right\}  $, and thus $0\leq k\leq n$.

Recall that $\left[  n\right]  =\left\{  1,2,\ldots,n\right\}  $ (by the
definition of $\left[  n\right]  $).

We have $\left(  i_{1},i_{2},\ldots,i_{k}\right)  \in\left[  1\right]
\times\left[  2\right]  \times\cdots\times\left[  k\right]  $. Thus,%
\begin{equation}
i_{j}\in\left[  j\right]  \ \ \ \ \ \ \ \ \ \ \text{for each }j\in\left\{
1,2,\ldots,k\right\}  . \label{pf.lem.sol.transpose.code.uni.1}%
\end{equation}

\textbf{(a)} Let $i\in\left\{  k+1,k+2,\ldots,n\right\}  $. Thus, $i\geq
k+1>k$.

We claim that%
\begin{equation}
\left(  t_{1,i_{1}}\circ t_{2,i_{2}}\circ\cdots\circ t_{p,i_{p}}\right)
\left(  i\right)  =i\ \ \ \ \ \ \ \ \ \ \text{for each }p\in\left\{
0,1,\ldots,k\right\}  . \label{pf.lem.sol.transpose.code.uni.c1.1}%
\end{equation}

[\textit{Proof of (\ref{pf.lem.sol.transpose.code.uni.c1.1}):} We shall prove
(\ref{pf.lem.sol.transpose.code.uni.c1.1}) by induction over $p$:

\textit{Induction base:} We have%
\[
\underbrace{\left(  t_{1,i_{1}}\circ t_{2,i_{2}}\circ\cdots\circ t_{0,i_{0}%
}\right)  }_{=\left(  \text{empty composition}\right)  =\operatorname*{id}%
}\left(  i\right)  =\operatorname*{id}\left(  i\right)  =i.
\]
In other words, (\ref{pf.lem.sol.transpose.code.uni.c1.1}) holds for $p=0$.
This completes the induction base.

\textit{Induction step:} Let $q\in\left\{  0,1,\ldots,k\right\}  $ be
positive. Assume that (\ref{pf.lem.sol.transpose.code.uni.c1.1}) holds for
$p=q-1$. We must prove that (\ref{pf.lem.sol.transpose.code.uni.c1.1}) holds
for $p=q$.

We have assumed that (\ref{pf.lem.sol.transpose.code.uni.c1.1}) holds for
$p=q-1$. In other words, we have $\left(  t_{1,i_{1}}\circ t_{2,i_{2}}%
\circ\cdots\circ t_{q-1,i_{q-1}}\right)  \left(  i\right)  =i$.

We have $q\in\left\{  0,1,\ldots,k\right\}  $, so that $q\leq k$. Hence,
$k\geq q$, so that $i>k\geq q$. Thus, $i\neq q$. Furthermore, $q\neq0$ (since
$q$ is positive). Combining $q\in\left\{  0,1,\ldots,k\right\}  $ with
$q\neq0$, we obtain $q\in\left\{  0,1,\ldots,k\right\}  \setminus\left\{
0\right\}  =\left\{  1,2,\ldots,k\right\}  $. Hence,
(\ref{pf.lem.sol.transpose.code.uni.1}) (applied to $j=q$) yields $i_{q}%
\in\left[  q\right]  =\left\{  1,2,\ldots,q\right\}  $ (by the definition of
$\left[  q\right]  $). Hence, $i_{q}\leq q$. Thus, $q\geq i_{q}$, so that
$i>q\geq i_{q}$. Thus, $i\neq i_{q}$. Finally, $i_{q}\in\left\{
1,2,\ldots,q\right\}  \subseteq\left\{  1,2,\ldots,n\right\}  $ (since $q\leq
k\leq n$). In other words, $i_{q}\in\left[  n\right]  $ (since $\left[
n\right]  =\left\{  1,2,\ldots,n\right\}  $).

Also, $q$ is positive. Thus, $q\in\left\{  1,2,\ldots,n\right\}  $ (since
$q\leq k\leq n$). In other words, $q\in\left[  n\right]  $ (since $\left[
n\right]  =\left\{  1,2,\ldots,n\right\}  $). Also,
\begin{align*}
i  &  \in\left\{  k+1,k+2,\ldots,n\right\}  \subseteq\left\{  1,2,\ldots
,n\right\}  \ \ \ \ \ \ \ \ \ \ \left(  \text{since }k\geq0\right) \\
&  =\left[  n\right]
\end{align*}
(since $\left[  n\right]  =\left\{  1,2,\ldots,n\right\}  $). Thus, $i$ is an
element of $\left[  n\right]  $ other than $q$ and $i_{q}$ (since $i\neq q$
and $i\neq i_{q}$). In other words, $i\in\left[  n\right]  \setminus\left\{
q,i_{q}\right\}  $. Hence, Lemma \ref{lem.sol.transpose.code.tij1}
\textbf{(c)} (applied to $q$, $i_{q}$ and $i$ instead of $i$, $j$ and $k$)
shows that $t_{q,i_{q}}\left(  i\right)  =i$. Now,%
\begin{align*}
&  \underbrace{\left(  t_{1,i_{1}}\circ t_{2,i_{2}}\circ\cdots\circ
t_{q,i_{q}}\right)  }_{=\left(  t_{1,i_{1}}\circ t_{2,i_{2}}\circ\cdots\circ
t_{q-1,i_{q-1}}\right)  \circ t_{q,i_{q}}}\left(  i\right) \\
&  =\left(  \left(  t_{1,i_{1}}\circ t_{2,i_{2}}\circ\cdots\circ
t_{q-1,i_{q-1}}\right)  \circ t_{q,i_{q}}\right)  \left(  i\right)  =\left(
t_{1,i_{1}}\circ t_{2,i_{2}}\circ\cdots\circ t_{q-1,i_{q-1}}\right)  \left(
\underbrace{t_{q,i_{q}}\left(  i\right)  }_{=i}\right) \\
&  =\left(  t_{1,i_{1}}\circ t_{2,i_{2}}\circ\cdots\circ t_{q-1,i_{q-1}%
}\right)  \left(  i\right)  =i.
\end{align*}
In other words, (\ref{pf.lem.sol.transpose.code.uni.c1.1}) holds for $p=q$.
This completes the induction step. Thus,
(\ref{pf.lem.sol.transpose.code.uni.c1.1}) is proven.]

Now, we can apply (\ref{pf.lem.sol.transpose.code.uni.c1.1}) to $p=k$ (since
$k\in\left\{  0,1,\ldots,k\right\}  $). We thus obtain%
\[
\left(  t_{1,i_{1}}\circ t_{2,i_{2}}\circ\cdots\circ t_{k,i_{k}}\right)
\left(  i\right)  =i.
\]
Hence,%
\[
\underbrace{\sigma}_{=t_{1,i_{1}}\circ t_{2,i_{2}}\circ\cdots\circ t_{k,i_{k}%
}}\left(  i\right)  =\left(  t_{1,i_{1}}\circ t_{2,i_{2}}\circ\cdots\circ
t_{k,i_{k}}\right)  \left(  i\right)  =i.
\]
This proves Lemma \ref{lem.sol.transpose.code.uni} \textbf{(a)}.

\textbf{(b)} Assume that $k>0$. Thus, $k$ is a positive integer. Hence,
$k-1\in\mathbb{N}$. Also, $k-1\leq k\leq n$, so that $k-1\in\left\{
0,1,\ldots,n\right\}  $ (since $k-1\in\mathbb{N}$). Furthermore, $k\in\left\{
\left(  k-1\right)  +1,\left(  k-1\right)  +2,\ldots,n\right\}  $ (because
$k=\left(  k-1\right)  +1\geq\left(  k-1\right)  +1$ and $k\leq n$).

Define $\tau\in S_{n}$ by $\tau=t_{1,i_{1}}\circ t_{2,i_{2}}\circ\cdots\circ
t_{k-1,i_{k-1}}$. (This is well-defined, since $k-1\in\mathbb{N}$.) Each
$j\in\left\{  1,2,\ldots,k-1\right\}  $ satisfies $j\in\left\{  1,2,\ldots
,k\right\}  $ (since $\left\{  1,2,\ldots,k-1\right\}  \subseteq\left\{
1,2,\ldots,k\right\}  $), and thus satisfies $i_{j}\in\left[  j\right]  $ (by
(\ref{pf.lem.sol.transpose.code.uni.1})). Hence, we have $\left(  i_{1}%
,i_{2},\ldots,i_{k-1}\right)  \in\left[  1\right]  \times\left[  2\right]
\times\cdots\times\left[  k-1\right]  $.

Hence, Lemma \ref{lem.sol.transpose.code.uni} \textbf{(a)} (applied to $k-1$,
$\tau$ and $k$ instead of $k$, $\sigma$ and $i$) yields $\tau\left(  k\right)
=k$.

We have $k>0$, so that $k\neq0$. Combining this with $k\in\left\{
0,1,\ldots,n\right\}  $, we find $k\in\left\{  0,1,\ldots,n\right\}
\setminus\left\{  0\right\}  =\left\{  1,2,\ldots,n\right\}  =\left[
n\right]  $. Also, $k\in\left\{  1,2,\ldots,k\right\}  $ (since $k>0$). Hence,
(\ref{pf.lem.sol.transpose.code.uni.1}) (applied to $j=k$) yields $i_{k}%
\in\left[  k\right]  =\left\{  1,2,\ldots,k\right\}  $ (by the definition of
$\left[  k\right]  $). But $k\leq n$ and thus $\left\{  1,2,\ldots,k\right\}
\subseteq\left\{  1,2,\ldots,n\right\}  =\left[  n\right]  $. Hence, $i_{k}%
\in\left\{  1,2,\ldots,k\right\}  \subseteq\left[  n\right]  $.

Thus, Lemma \ref{lem.sol.transpose.code.tij1} \textbf{(b)} (applied to $k$ and
$i_{k}$ instead of $i$ and $j$) shows that $t_{k,i_{k}}\left(  i_{k}\right)
=k$. Now,%
\begin{align*}
\sigma &  =t_{1,i_{1}}\circ t_{2,i_{2}}\circ\cdots\circ t_{k,i_{k}%
}=\underbrace{\left(  t_{1,i_{1}}\circ t_{2,i_{2}}\circ\cdots\circ
t_{k-1,i_{k-1}}\right)  }_{=\tau}\circ t_{k,i_{k}}\ \ \ \ \ \ \ \ \ \ \left(
\text{since }k>0\right) \\
&  =\tau\circ t_{k,i_{k}}.
\end{align*}
Applying both sides of this equality to $i_{k}$, we find%
\[
\sigma\left(  i_{k}\right)  =\left(  \tau\circ t_{k,i_{k}}\right)  \left(
i_{k}\right)  =\tau\left(  \underbrace{t_{k,i_{k}}\left(  i_{k}\right)  }%
_{=k}\right)  =\tau\left(  k\right)  =k.
\]
This proves Lemma \ref{lem.sol.transpose.code.uni} \textbf{(b)}.
\end{proof}
\end{verlong}

\begin{lemma}
\label{lem.sol.transpose.code.uni2}Let $n\in\mathbb{N}$. Let $k\in\left\{
0,1,\ldots,n\right\}  $. Let $\left(  u_{1},u_{2},\ldots,u_{k}\right)
\in\left[  1\right]  \times\left[  2\right]  \times\cdots\times\left[
k\right]  $ and $\left(  v_{1},v_{2},\ldots,v_{k}\right)  \in\left[  1\right]
\times\left[  2\right]  \times\cdots\times\left[  k\right]  $ be two
$k$-tuples satisfying
\[
t_{1,u_{1}}\circ t_{2,u_{2}}\circ\cdots\circ t_{k,u_{k}}=t_{1,v_{1}}\circ
t_{2,v_{2}}\circ\cdots\circ t_{k,v_{k}}.
\]
Then, $\left(  u_{1},u_{2},\ldots,u_{k}\right)  =\left(  v_{1},v_{2}%
,\ldots,v_{k}\right)  $.
\end{lemma}

\begin{proof}
[Proof of Lemma \ref{lem.sol.transpose.code.uni2}.]We shall prove Lemma
\ref{lem.sol.transpose.code.uni2} by induction over $k$:

\begin{vershort}
\textit{Induction base:} Lemma \ref{lem.sol.transpose.code.uni2} is true when
$k=0$\ \ \ \ \footnote{\textit{Proof.} There exists only one $0$-tuple:
namely, the empty $0$-tuple $\left(  {}\right)  $. Thus, any two $0$-tuples
are equal.
\par
Lemma \ref{lem.sol.transpose.code.uni2} claims the equality of two $k$-tuples.
When $k=0$, any two $k$-tuples are equal (since any two $0$-tuples are equal).
Hence, Lemma \ref{lem.sol.transpose.code.uni2} is true when $k=0$.}.
\end{vershort}

\begin{verlong}
\textit{Induction base:} Lemma \ref{lem.sol.transpose.code.uni2} is true when
$k=0$\ \ \ \ \footnote{\textit{Proof.} Let $n$, $k$, $\left(  u_{1}%
,u_{2},\ldots,u_{k}\right)  $ and $\left(  v_{1},v_{2},\ldots,v_{k}\right)  $
be as in Lemma \ref{lem.sol.transpose.code.uni2}. Assume that $k=0$. We have
to prove that Lemma \ref{lem.sol.transpose.code.uni2} holds under this
assumption.
\par
Recall that there exists only one $0$-tuple: namely, the empty $0$-tuple
$\left(  {}\right)  $.
\par
From $k=0$, we obtain $\left(  u_{1},u_{2},\ldots,u_{k}\right)  =\left(
u_{1},u_{2},\ldots,u_{0}\right)  =\left(  {}\right)  $ (this is the empty
$0$-tuple). Also, from $k=0$, we obtain $\left(  v_{1},v_{2},\ldots
,v_{k}\right)  =\left(  v_{1},v_{2},\ldots,v_{0}\right)  =\left(  {}\right)
$. Thus, $\left(  u_{1},u_{2},\ldots,u_{k}\right)  =\left(  {}\right)
=\left(  v_{1},v_{2},\ldots,v_{k}\right)  $.
\par
Now, forget that we fixed $n$, $k$, $\left(  u_{1},u_{2},\ldots,u_{k}\right)
$ and $\left(  v_{1},v_{2},\ldots,v_{k}\right)  $. We thus have proven that if
$n$, $k$, $\left(  u_{1},u_{2},\ldots,u_{k}\right)  $ and $\left(  v_{1}%
,v_{2},\ldots,v_{k}\right)  $ are as in Lemma
\ref{lem.sol.transpose.code.uni2}, and if we have $k=0$, then $\left(
u_{1},u_{2},\ldots,u_{k}\right)  =\left(  v_{1},v_{2},\ldots,v_{k}\right)  $.
In other words, Lemma \ref{lem.sol.transpose.code.uni2} is true when $k=0$.
Qed.}. This completes the induction base.
\end{verlong}

\textit{Induction step:} Let $K\in\left\{  0,1,\ldots,n\right\}  $ be
positive. Assume that Lemma \ref{lem.sol.transpose.code.uni2} holds when
$k=K-1$. We must show that Lemma \ref{lem.sol.transpose.code.uni2} holds when
$k=K$.

We have assumed that Lemma \ref{lem.sol.transpose.code.uni2} holds when
$k=K-1$. In other words, the following statement holds:

\begin{statement}
\textit{Statement 1:} Let $n\in\mathbb{N}$. Assume that $K-1\in\left\{
0,1,\ldots,n\right\}  $. Let \newline$\left(  u_{1},u_{2},\ldots
,u_{K-1}\right)  \in\left[  1\right]  \times\left[  2\right]  \times
\cdots\times\left[  K-1\right]  $ and $\left(  v_{1},v_{2},\ldots
,v_{K-1}\right)  \in\left[  1\right]  \times\left[  2\right]  \times
\cdots\times\left[  K-1\right]  $ be two $\left(  K-1\right)  $-tuples
satisfying
\[
t_{1,u_{1}}\circ t_{2,u_{2}}\circ\cdots\circ t_{K-1,u_{K-1}}=t_{1,v_{1}}\circ
t_{2,v_{2}}\circ\cdots\circ t_{K-1,v_{K-1}}.
\]
Then, $\left(  u_{1},u_{2},\ldots,u_{K-1}\right)  =\left(  v_{1},v_{2}%
,\ldots,v_{K-1}\right)  $.
\end{statement}

We must show that Lemma \ref{lem.sol.transpose.code.uni2} holds when $k=K$. In
other words, we must prove the following statement:

\begin{statement}
\textit{Statement 2:} Let $n\in\mathbb{N}$. Assume that $K\in\left\{
0,1,\ldots,n\right\}  $. Let $\left(  u_{1},u_{2},\ldots,u_{K}\right)
\in\left[  1\right]  \times\left[  2\right]  \times\cdots\times\left[
K\right]  $ and $\left(  v_{1},v_{2},\ldots,v_{K}\right)  \in\left[  1\right]
\times\left[  2\right]  \times\cdots\times\left[  K\right]  $ be two
$K$-tuples satisfying
\[
t_{1,u_{1}}\circ t_{2,u_{2}}\circ\cdots\circ t_{K,u_{K}}=t_{1,v_{1}}\circ
t_{2,v_{2}}\circ\cdots\circ t_{K,v_{K}}.
\]
Then, $\left(  u_{1},u_{2},\ldots,u_{K}\right)  =\left(  v_{1},v_{2}%
,\ldots,v_{K}\right)  $.
\end{statement}

\begin{vershort}
[\textit{Proof of Statement 2:} We have $K\in\left\{  0,1,\ldots,n\right\}  $.
Thus, $K\in\left\{  1,2,\ldots,n\right\}  $ (since $K$ is positive), so that
$K-1\in\left\{  0,1,\ldots,n-1\right\}  \subseteq\left\{  0,1,\ldots
,n\right\}  $.

We have $\left(  u_{1},u_{2},\ldots,u_{K}\right)  \in\left[  1\right]
\times\left[  2\right]  \times\cdots\times\left[  K\right]  $. In other words,
$u_{j}\in\left[  j\right]  $ for each $j\in\left\{  1,2,\ldots,K\right\}  $.
Hence, $\left(  u_{1},u_{2},\ldots,u_{K-1}\right)  \in\left[  1\right]
\times\left[  2\right]  \times\cdots\times\left[  K-1\right]  $. The same
argument (applied to $v_{i}$ instead of $u_{i}$) shows that $\left(
v_{1},v_{2},\ldots,v_{K-1}\right)  \in\left[  1\right]  \times\left[
2\right]  \times\cdots\times\left[  K-1\right]  $.

Define $\sigma\in S_{n}$ by $\sigma=t_{1,u_{1}}\circ t_{2,u_{2}}\circ
\cdots\circ t_{K,u_{K}}$. Hence, Lemma \ref{lem.sol.transpose.code.uni}
\textbf{(b)} (applied to $k=K$ and $i_{j}=u_{j}$) yields that $\sigma\left(
u_{K}\right)  =K$ (since $K>0$). But we also have $\sigma=t_{1,u_{1}}\circ
t_{2,u_{2}}\circ\cdots\circ t_{K,u_{K}}=t_{1,v_{1}}\circ t_{2,v_{2}}%
\circ\cdots\circ t_{K,v_{K}}$. Hence, Lemma \ref{lem.sol.transpose.code.uni}
\textbf{(b)} (applied to $k=K$ and $i_{j}=v_{j}$) yields that $\sigma\left(
v_{K}\right)  =K$ (since $K>0$). Hence, $\sigma\left(  u_{K}\right)
=K=\sigma\left(  v_{K}\right)  $.

We have $\sigma\in S_{n}$. In other words, $\sigma$ is a permutation of
$\left[  n\right]  $ (since $S_{n}$ is the set of all permutations of $\left[
n\right]  $). In other words, $\sigma$ is a bijection $\left[  n\right]
\rightarrow\left[  n\right]  $. Hence, this map $\sigma$ is bijective, thus
injective. Therefore, from $\sigma\left(  u_{K}\right)  =\sigma\left(
v_{K}\right)  $, we conclude that $u_{K}=v_{K}$.

Recall that $u_{j}\in\left[  j\right]  $ for each $j\in\left\{  1,2,\ldots
,K\right\}  $. Applying this to $j=K$, we obtain $u_{K}\in\left[  K\right]  $.
In view of $u_{K}=v_{K}$, this rewrites as $v_{K}\in\left[  K\right]  $. But
$K\leq n$ (since $K\in\left\{  0,1,\ldots,n\right\}  $), so that $\left[
K\right]  \subseteq\left[  n\right]  $. Hence, $v_{K}\in\left[  K\right]
\subseteq\left[  n\right]  $. Furthermore, $K\in\left\{  1,2,\ldots,n\right\}
=\left[  n\right]  $. Hence, Lemma \ref{lem.sol.transpose.code.tij1}
\textbf{(d)} (applied to $K$ and $v_{K}$ instead of $i$ and $j$) yields
$t_{K,v_{K}}\circ t_{K,v_{K}}=\operatorname*{id}$.

Now, recall that $K>0$. Hence,%
\begin{align*}
&  t_{1,u_{1}}\circ t_{2,u_{2}}\circ\cdots\circ t_{K,u_{K}}\\
&  =\left(  t_{1,u_{1}}\circ t_{2,u_{2}}\circ\cdots\circ t_{K-1,u_{K-1}%
}\right)  \circ\underbrace{t_{K,u_{K}}}_{\substack{=t_{K,v_{K}}\\\text{(since
}u_{K}=v_{K}\text{)}}}=\left(  t_{1,u_{1}}\circ t_{2,u_{2}}\circ\cdots\circ
t_{K-1,u_{K-1}}\right)  \circ t_{K,v_{K}}.
\end{align*}
Comparing this with%
\begin{align*}
&  t_{1,u_{1}}\circ t_{2,u_{2}}\circ\cdots\circ t_{K,u_{K}}\\
&  =t_{1,v_{1}}\circ t_{2,v_{2}}\circ\cdots\circ t_{K,v_{K}}=\left(
t_{1,v_{1}}\circ t_{2,v_{2}}\circ\cdots\circ t_{K-1,v_{K-1}}\right)  \circ
t_{K,v_{K}},
\end{align*}
we obtain%
\[
\left(  t_{1,u_{1}}\circ t_{2,u_{2}}\circ\cdots\circ t_{K-1,u_{K-1}}\right)
\circ t_{K,v_{K}}=\left(  t_{1,v_{1}}\circ t_{2,v_{2}}\circ\cdots\circ
t_{K-1,v_{K-1}}\right)  \circ t_{K,v_{K}}.
\]
Thus,%
\begin{align*}
&  \underbrace{\left(  t_{1,u_{1}}\circ t_{2,u_{2}}\circ\cdots\circ
t_{K-1,u_{K-1}}\right)  \circ t_{K,v_{K}}}_{=\left(  t_{1,v_{1}}\circ
t_{2,v_{2}}\circ\cdots\circ t_{K-1,v_{K-1}}\right)  \circ t_{K,v_{K}}}\circ
t_{K,v_{K}}\\
&  =\left(  t_{1,v_{1}}\circ t_{2,v_{2}}\circ\cdots\circ t_{K-1,v_{K-1}%
}\right)  \circ\underbrace{t_{K,v_{K}}\circ t_{K,v_{K}}}_{=\operatorname*{id}%
}\\
&  =t_{1,v_{1}}\circ t_{2,v_{2}}\circ\cdots\circ t_{K-1,v_{K-1}}.
\end{align*}
Comparing this with%
\[
\left(  t_{1,u_{1}}\circ t_{2,u_{2}}\circ\cdots\circ t_{K-1,u_{K-1}}\right)
\circ\underbrace{t_{K,v_{K}}\circ t_{K,v_{K}}}_{=\operatorname*{id}%
}=t_{1,u_{1}}\circ t_{2,u_{2}}\circ\cdots\circ t_{K-1,u_{K-1}},
\]
we obtain%
\[
t_{1,u_{1}}\circ t_{2,u_{2}}\circ\cdots\circ t_{K-1,u_{K-1}}=t_{1,v_{1}}\circ
t_{2,v_{2}}\circ\cdots\circ t_{K-1,v_{K-1}}.
\]
Hence, Statement 1 shows that $\left(  u_{1},u_{2},\ldots,u_{K-1}\right)
=\left(  v_{1},v_{2},\ldots,v_{K-1}\right)  $. In other words,%
\[
u_{j}=v_{j}\ \ \ \ \ \ \ \ \ \ \text{for each }j\in\left\{  1,2,\ldots
,K-1\right\}  .
\]
Combining this with $u_{K}=v_{K}$, we conclude that
\[
u_{j}=v_{j}\ \ \ \ \ \ \ \ \ \ \text{for each }j\in\left\{  1,2,\ldots
,K\right\}  .
\]
In other words, $\left(  u_{1},u_{2},\ldots,u_{K}\right)  =\left(  v_{1}%
,v_{2},\ldots,v_{K}\right)  $. This proves Statement 2.]
\end{vershort}

\begin{verlong}
[\textit{Proof of Statement 2:} We have $\left[  n\right]  =\left\{
1,2,\ldots,n\right\}  $ (by the definition of $\left[  n\right]  $).

We have $K>0$ (since $K$ is positive), so that $K\neq0$. Combining
$K\in\left\{  0,1,\ldots,n\right\}  $ with $K\neq0$, we obtain $K\in\left\{
0,1,\ldots,n\right\}  \setminus\left\{  0\right\}  =\left\{  1,2,\ldots
,n\right\}  $, so that $K-1\in\left\{  0,1,\ldots,n-1\right\}  \subseteq
\left\{  0,1,\ldots,n\right\}  $.

We have $\left(  u_{1},u_{2},\ldots,u_{K}\right)  \in\left[  1\right]
\times\left[  2\right]  \times\cdots\times\left[  K\right]  $. Therefore,
$\left(  u_{1},u_{2},\ldots,u_{K-1}\right)  \in\left[  1\right]  \times\left[
2\right]  \times\cdots\times\left[  K-1\right]  $%
\ \ \ \ \footnote{\textit{Proof.} Each $j\in\left\{  1,2,\ldots,K-1\right\}  $
satisfies $j\in\left\{  1,2,\ldots,K\right\}  $ (since $\left\{
1,2,\ldots,K-1\right\}  \subseteq\left\{  1,2,\ldots,K\right\}  $).
\par
We have $\left(  u_{1},u_{2},\ldots,u_{K}\right)  \in\left[  1\right]
\times\left[  2\right]  \times\cdots\times\left[  K\right]  $. In other words,
each $j\in\left\{  1,2,\ldots,K\right\}  $ satisfies $u_{j}\in\left[
j\right]  $. Hence, each $j\in\left\{  1,2,\ldots,K-1\right\}  $ satisfies
$u_{j}\in\left[  j\right]  $ (since each $j\in\left\{  1,2,\ldots,K-1\right\}
$ satisfies $j\in\left\{  1,2,\ldots,K\right\}  $). In other words, $\left(
u_{1},u_{2},\ldots,u_{K-1}\right)  \in\left[  1\right]  \times\left[
2\right]  \times\cdots\times\left[  K-1\right]  $. Qed.}. The same argument
(applied to $v_{i}$ instead of $u_{i}$) shows that $\left(  v_{1},v_{2}%
,\ldots,v_{K-1}\right)  \in\left[  1\right]  \times\left[  2\right]
\times\cdots\times\left[  K-1\right]  $.

Define $\sigma\in S_{n}$ by $\sigma=t_{1,u_{1}}\circ t_{2,u_{2}}\circ
\cdots\circ t_{K,u_{K}}$. Hence, Lemma \ref{lem.sol.transpose.code.uni}
\textbf{(b)} (applied to $k=K$ and $i_{j}=u_{j}$) yields that $\sigma\left(
u_{K}\right)  =K$ (since $K>0$). But we also have $\sigma=t_{1,u_{1}}\circ
t_{2,u_{2}}\circ\cdots\circ t_{K,u_{K}}=t_{1,v_{1}}\circ t_{2,v_{2}}%
\circ\cdots\circ t_{K,v_{K}}$. Hence, Lemma \ref{lem.sol.transpose.code.uni}
\textbf{(b)} (applied to $k=K$ and $i_{j}=v_{j}$) yields that $\sigma\left(
v_{K}\right)  =K$ (since $K>0$). Hence, $\sigma\left(  u_{K}\right)
=K=\sigma\left(  v_{K}\right)  $.

We have $\sigma\in S_{n}$. In other words, $\sigma$ is a permutation of
$\left\{  1,2,\ldots,n\right\}  $ (since $S_{n}$ is the set of all
permutations of $\left\{  1,2,\ldots,n\right\}  $). In other words, $\sigma$
is a permutation of $\left[  n\right]  $ (since $\left[  n\right]  =\left\{
1,2,\ldots,n\right\}  $). In other words, $\sigma$ is a bijection $\left[
n\right]  \rightarrow\left[  n\right]  $. Hence, this map $\sigma$ is
bijective, thus injective. Therefore, from $\sigma\left(  u_{K}\right)
=\sigma\left(  v_{K}\right)  $, we conclude that $u_{K}=v_{K}$.

We have $K\in\left\{  1,2,\ldots,K\right\}  $ (since $K$ is a positive
integer). But $\left(  v_{1},v_{2},\ldots,v_{K}\right)  \in\left[  1\right]
\times\left[  2\right]  \times\cdots\times\left[  K\right]  $. In other words,
$v_{j}\in\left[  j\right]  $ for each $j\in\left\{  1,2,\ldots,K\right\}  $.
Applying this to $j=K$, we obtain $v_{K}\in\left[  K\right]  $ (since
$K\in\left\{  1,2,\ldots,K\right\}  $). But $K\leq n$ (since $K\in\left\{
0,1,\ldots,n\right\}  $), so that $\left\{  1,2,\ldots,K\right\}
\subseteq\left\{  1,2,\ldots,n\right\}  =\left[  n\right]  $. The definition
of $\left[  K\right]  $ yields $\left[  K\right]  =\left\{  1,2,\ldots
,K\right\}  \subseteq\left[  n\right]  $. Hence, $v_{K}\in\left[  K\right]
\subseteq\left[  n\right]  $. Furthermore, $K\in\left\{  1,2,\ldots,n\right\}
=\left[  n\right]  $. Hence, Lemma \ref{lem.sol.transpose.code.tij1}
\textbf{(d)} (applied to $K$ and $v_{K}$ instead of $i$ and $j$) yields
$t_{K,v_{K}}\circ t_{K,v_{K}}=\operatorname*{id}$.

Now, recall that $K>0$. Hence,%
\begin{align*}
t_{1,u_{1}}\circ t_{2,u_{2}}\circ\cdots\circ t_{K,u_{K}}  &  =\left(
t_{1,u_{1}}\circ t_{2,u_{2}}\circ\cdots\circ t_{K-1,u_{K-1}}\right)
\circ\underbrace{t_{K,u_{K}}}_{\substack{=t_{K,v_{K}}\\\text{(since }%
u_{K}=v_{K}\text{)}}}\\
&  =\left(  t_{1,u_{1}}\circ t_{2,u_{2}}\circ\cdots\circ t_{K-1,u_{K-1}%
}\right)  \circ t_{K,v_{K}}.
\end{align*}
Comparing this with%
\begin{align*}
t_{1,u_{1}}\circ t_{2,u_{2}}\circ\cdots\circ t_{K,u_{K}}  &  =t_{1,v_{1}}\circ
t_{2,v_{2}}\circ\cdots\circ t_{K,v_{K}}\\
&  =\left(  t_{1,v_{1}}\circ t_{2,v_{2}}\circ\cdots\circ t_{K-1,v_{K-1}%
}\right)  \circ t_{K,v_{K}}\ \ \ \ \ \ \ \ \ \ \left(  \text{since
}K>0\right)  ,
\end{align*}
we obtain%
\[
\left(  t_{1,u_{1}}\circ t_{2,u_{2}}\circ\cdots\circ t_{K-1,u_{K-1}}\right)
\circ t_{K,v_{K}}=\left(  t_{1,v_{1}}\circ t_{2,v_{2}}\circ\cdots\circ
t_{K-1,v_{K-1}}\right)  \circ t_{K,v_{K}}.
\]
Thus,%
\begin{align*}
&  \underbrace{\left(  t_{1,u_{1}}\circ t_{2,u_{2}}\circ\cdots\circ
t_{K-1,u_{K-1}}\right)  \circ t_{K,v_{K}}}_{=\left(  t_{1,v_{1}}\circ
t_{2,v_{2}}\circ\cdots\circ t_{K-1,v_{K-1}}\right)  \circ t_{K,v_{K}}}\circ
t_{K,v_{K}}\\
&  =\left(  t_{1,v_{1}}\circ t_{2,v_{2}}\circ\cdots\circ t_{K-1,v_{K-1}%
}\right)  \circ\underbrace{t_{K,v_{K}}\circ t_{K,v_{K}}}_{=\operatorname*{id}%
}\\
&  =t_{1,v_{1}}\circ t_{2,v_{2}}\circ\cdots\circ t_{K-1,v_{K-1}}.
\end{align*}
Comparing this with%
\[
\left(  t_{1,u_{1}}\circ t_{2,u_{2}}\circ\cdots\circ t_{K-1,u_{K-1}}\right)
\circ\underbrace{t_{K,v_{K}}\circ t_{K,v_{K}}}_{=\operatorname*{id}%
}=t_{1,u_{1}}\circ t_{2,u_{2}}\circ\cdots\circ t_{K-1,u_{K-1}},
\]
we obtain%
\[
t_{1,u_{1}}\circ t_{2,u_{2}}\circ\cdots\circ t_{K-1,u_{K-1}}=t_{1,v_{1}}\circ
t_{2,v_{2}}\circ\cdots\circ t_{K-1,v_{K-1}}.
\]
Hence, Statement 1 shows that $\left(  u_{1},u_{2},\ldots,u_{K-1}\right)
=\left(  v_{1},v_{2},\ldots,v_{K-1}\right)  $. In other words,%
\begin{equation}
u_{j}=v_{j}\ \ \ \ \ \ \ \ \ \ \text{for each }j\in\left\{  1,2,\ldots
,K-1\right\}  . \label{pf.lem.sol.transpose.code.uni2.uj=vj}%
\end{equation}

Thus,%
\begin{equation}
u_{j}=v_{j}\ \ \ \ \ \ \ \ \ \ \text{for each }j\in\left\{  1,2,\ldots
,K\right\}  \label{pf.lem.sol.transpose.code.uni2.uj=vj2}%
\end{equation}
\footnote{\textit{Proof of (\ref{pf.lem.sol.transpose.code.uni2.uj=vj2}):} Let
$j\in\left\{  1,2,\ldots,K\right\}  $. We must prove that $u_{j}=v_{j}$.
\par
If $j\in\left\{  1,2,\ldots,K-1\right\}  $, then this follows immediately from
(\ref{pf.lem.sol.transpose.code.uni2.uj=vj}). Thus, for the rest of this
proof, we can WLOG assume that we don't have $j\in\left\{  1,2,\ldots
,K-1\right\}  $. Assume this.
\par
We have $j\notin\left\{  1,2,\ldots,K-1\right\}  $ (since we don't have
$j\in\left\{  1,2,\ldots,K-1\right\}  $). Combining this with $j\in\left\{
1,2,\ldots,K\right\}  $, we obtain $j\in\left\{  1,2,\ldots,K\right\}
\setminus\left\{  1,2,\ldots,K-1\right\}  \subseteq\left\{  K\right\}  $. In
other words, $j=K$. Now, recall that $u_{K}=v_{K}$. In view of $j=K$, this
rewrites as $u_{j}=v_{j}$. Thus, we have proven that $u_{j}=v_{j}$. This
completes the proof of (\ref{pf.lem.sol.transpose.code.uni2.uj=vj2}).}. In
other words, $\left(  u_{1},u_{2},\ldots,u_{K}\right)  =\left(  v_{1}%
,v_{2},\ldots,v_{K}\right)  $. This proves Statement 2.]
\end{verlong}

So Statement 2 is proven. In other words, Lemma
\ref{lem.sol.transpose.code.uni2} holds when $k=K$ (since Statement 2 is
precisely the claim of Lemma \ref{lem.sol.transpose.code.uni2} for $k=K$).
This completes the induction step. Thus, Lemma
\ref{lem.sol.transpose.code.uni2} is proven.
\end{proof}

\subsubsection{Solving Exercise \ref{exe.transpos.code}}

\begin{vershort}
\begin{proof}
[Solution to Exercise \ref{exe.transpos.code}.]We have $n\in\left\{
0,1,\ldots,n\right\}  $ (since $n\in\mathbb{N}$).

We have $\left(  \sigma\left(  i\right)  =i\text{ for each }i\in\left\{
n+1,n+2,\ldots,n\right\}  \right)  $. (Indeed, this is vacuously true, because
there exists no $i\in\left\{  n+1,n+2,\ldots,n\right\}  $.) Hence, Lemma
\ref{lem.sol.transpose.code.exist} (applied to $k=n$) yields that there exists
a $n$-tuple $\left(  i_{1},i_{2},\ldots,i_{n}\right)  \in\left[  1\right]
\times\left[  2\right]  \times\cdots\times\left[  n\right]  $ such that
$\sigma=t_{1,i_{1}}\circ t_{2,i_{2}}\circ\cdots\circ t_{n,i_{n}}$.

Moreover, there exists at most one such $n$-tuple\footnote{\textit{Proof.} Let
$\left(  u_{1},u_{2},\ldots,u_{n}\right)  $ and $\left(  v_{1},v_{2}%
,\ldots,v_{n}\right)  $ be two such $n$-tuples. We shall prove that $\left(
u_{1},u_{2},\ldots,u_{n}\right)  =\left(  v_{1},v_{2},\ldots,v_{n}\right)  $.
\par
We know that $\left(  u_{1},u_{2},\ldots,u_{n}\right)  $ is an $n$-tuple
$\left(  i_{1},i_{2},\ldots,i_{n}\right)  \in\left[  1\right]  \times\left[
2\right]  \times\cdots\times\left[  n\right]  $ such that $\sigma=t_{1,i_{1}%
}\circ t_{2,i_{2}}\circ\cdots\circ t_{n,i_{n}}$. In other words, $\left(
u_{1},u_{2},\ldots,u_{n}\right)  $ is an $n$-tuple in $\left[  1\right]
\times\left[  2\right]  \times\cdots\times\left[  n\right]  $ and satisfies
$\sigma=t_{1,u_{1}}\circ t_{2,u_{2}}\circ\cdots\circ t_{n,u_{n}}$. Thus,
$t_{1,u_{1}}\circ t_{2,u_{2}}\circ\cdots\circ t_{n,u_{n}}=\sigma$ and $\left(
u_{1},u_{2},\ldots,u_{n}\right)  \in\left[  1\right]  \times\left[  2\right]
\times\cdots\times\left[  n\right]  $ (since $\left(  u_{1},u_{2},\ldots
,u_{n}\right)  $ is an $n$-tuple in $\left[  1\right]  \times\left[  2\right]
\times\cdots\times\left[  n\right]  $).
\par
We know that $\left(  v_{1},v_{2},\ldots,v_{n}\right)  $ is an $n$-tuple
$\left(  i_{1},i_{2},\ldots,i_{n}\right)  \in\left[  1\right]  \times\left[
2\right]  \times\cdots\times\left[  n\right]  $ such that $\sigma=t_{1,i_{1}%
}\circ t_{2,i_{2}}\circ\cdots\circ t_{n,i_{n}}$. In other words, $\left(
v_{1},v_{2},\ldots,v_{n}\right)  $ is an $n$-tuple in $\left[  1\right]
\times\left[  2\right]  \times\cdots\times\left[  n\right]  $ and satisfies
$\sigma=t_{1,v_{1}}\circ t_{2,v_{2}}\circ\cdots\circ t_{n,v_{n}}$. Thus,
$\left(  v_{1},v_{2},\ldots,v_{n}\right)  \in\left[  1\right]  \times\left[
2\right]  \times\cdots\times\left[  n\right]  $ (since $\left(  v_{1}%
,v_{2},\ldots,v_{n}\right)  $ is an $n$-tuple in $\left[  1\right]
\times\left[  2\right]  \times\cdots\times\left[  n\right]  $).
\par
Now, $t_{1,u_{1}}\circ t_{2,u_{2}}\circ\cdots\circ t_{n,u_{n}}=\sigma
=t_{1,v_{1}}\circ t_{2,v_{2}}\circ\cdots\circ t_{n,v_{n}}$. Hence, Lemma
\ref{lem.sol.transpose.code.uni2} (applied to $k=n$) yields $\left(
u_{1},u_{2},\ldots,u_{n}\right)  =\left(  v_{1},v_{2},\ldots,v_{n}\right)  $
(since $n\in\left\{  0,1,\ldots,n\right\}  $).
\par
Now, forget that we fixed $\left(  u_{1},u_{2},\ldots,u_{n}\right)  $ and
$\left(  v_{1},v_{2},\ldots,v_{n}\right)  $. We thus have proven that if
$\left(  u_{1},u_{2},\ldots,u_{n}\right)  $ and $\left(  v_{1},v_{2}%
,\ldots,v_{n}\right)  $ are two $n$-tuples $\left(  i_{1},i_{2},\ldots
,i_{n}\right)  \in\left[  1\right]  \times\left[  2\right]  \times\cdots
\times\left[  n\right]  $ such that $\sigma=t_{1,i_{1}}\circ t_{2,i_{2}}%
\circ\cdots\circ t_{n,i_{n}}$, then $\left(  u_{1},u_{2},\ldots,u_{n}\right)
=\left(  v_{1},v_{2},\ldots,v_{n}\right)  $. In other words, there exists at
most one $n$-tuple $\left(  i_{1},i_{2},\ldots,i_{n}\right)  \in\left[
1\right]  \times\left[  2\right]  \times\cdots\times\left[  n\right]  $ such
that $\sigma=t_{1,i_{1}}\circ t_{2,i_{2}}\circ\cdots\circ t_{n,i_{n}}$. Qed.}.
Hence, there is a unique $n$-tuple $\left(  i_{1},i_{2},\ldots,i_{n}\right)
\in\left[  1\right]  \times\left[  2\right]  \times\cdots\times\left[
n\right]  $ such that $\sigma=t_{1,i_{1}}\circ t_{2,i_{2}}\circ\cdots\circ
t_{n,i_{n}}$ (because there exists such an $n$-tuple, and there exists at most
one such $n$-tuple). This solves Exercise \ref{exe.transpos.code}.
\end{proof}
\end{vershort}

\begin{verlong}
\begin{proof}
[Solution to Exercise \ref{exe.transpos.code}.]We have $n\in\left\{
0,1,\ldots,n\right\}  $ (since $n\in\mathbb{N}$).

Now, we make the following two claims:

\begin{statement}
\textit{Claim 1:} There exists \textbf{at least one} $n$-tuple $\left(
i_{1},i_{2},\ldots,i_{n}\right)  \in\left[  1\right]  \times\left[  2\right]
\times\cdots\times\left[  n\right]  $ such that $\sigma=t_{1,i_{1}}\circ
t_{2,i_{2}}\circ\cdots\circ t_{n,i_{n}}$.
\end{statement}

\begin{statement}
\textit{Claim 2:} There exists \textbf{at most one} $n$-tuple $\left(
i_{1},i_{2},\ldots,i_{n}\right)  \in\left[  1\right]  \times\left[  2\right]
\times\cdots\times\left[  n\right]  $ such that $\sigma=t_{1,i_{1}}\circ
t_{2,i_{2}}\circ\cdots\circ t_{n,i_{n}}$.
\end{statement}

[\textit{Proof of Claim 1:} We have $\left\{  n+1,n+2,\ldots,n\right\}
=\varnothing$. In other words, the set $\left\{  n+1,n+2,\ldots,n\right\}  $
is empty. In other words, there exists no $i\in\left\{  n+1,n+2,\ldots
,n\right\}  $.

We have $\left(  \sigma\left(  i\right)  =i\text{ for each }i\in\left\{
n+1,n+2,\ldots,n\right\}  \right)  $. (Indeed, this is vacuously true, because
there exists no $i\in\left\{  n+1,n+2,\ldots,n\right\}  $.) Hence, Lemma
\ref{lem.sol.transpose.code.exist} (applied to $k=n$) yields that there exists
a $n$-tuple $\left(  i_{1},i_{2},\ldots,i_{n}\right)  \in\left[  1\right]
\times\left[  2\right]  \times\cdots\times\left[  n\right]  $ such that
$\sigma=t_{1,i_{1}}\circ t_{2,i_{2}}\circ\cdots\circ t_{n,i_{n}}$. In other
words, there exists \textbf{at least one} $n$-tuple $\left(  i_{1}%
,i_{2},\ldots,i_{n}\right)  \in\left[  1\right]  \times\left[  2\right]
\times\cdots\times\left[  n\right]  $ such that $\sigma=t_{1,i_{1}}\circ
t_{2,i_{2}}\circ\cdots\circ t_{n,i_{n}}$. This proves Claim 1.]

[\textit{Proof of Claim 2:} Let $\left(  u_{1},u_{2},\ldots,u_{n}\right)  $
and $\left(  v_{1},v_{2},\ldots,v_{n}\right)  $ be two $n$-tuples $\left(
i_{1},i_{2},\ldots,i_{n}\right)  \in\left[  1\right]  \times\left[  2\right]
\times\cdots\times\left[  n\right]  $ such that $\sigma=t_{1,i_{1}}\circ
t_{2,i_{2}}\circ\cdots\circ t_{n,i_{n}}$. We shall prove that $\left(
u_{1},u_{2},\ldots,u_{n}\right)  =\left(  v_{1},v_{2},\ldots,v_{n}\right)  $.

We know that $\left(  u_{1},u_{2},\ldots,u_{n}\right)  $ is an $n$-tuple
$\left(  i_{1},i_{2},\ldots,i_{n}\right)  \in\left[  1\right]  \times\left[
2\right]  \times\cdots\times\left[  n\right]  $ such that $\sigma=t_{1,i_{1}%
}\circ t_{2,i_{2}}\circ\cdots\circ t_{n,i_{n}}$. In other words, $\left(
u_{1},u_{2},\ldots,u_{n}\right)  $ is an $n$-tuple in $\left[  1\right]
\times\left[  2\right]  \times\cdots\times\left[  n\right]  $ and satisfies
$\sigma=t_{1,u_{1}}\circ t_{2,u_{2}}\circ\cdots\circ t_{n,u_{n}}$. Thus,
$t_{1,u_{1}}\circ t_{2,u_{2}}\circ\cdots\circ t_{n,u_{n}}=\sigma$ and $\left(
u_{1},u_{2},\ldots,u_{n}\right)  \in\left[  1\right]  \times\left[  2\right]
\times\cdots\times\left[  n\right]  $ (since $\left(  u_{1},u_{2},\ldots
,u_{n}\right)  $ is an $n$-tuple in $\left[  1\right]  \times\left[  2\right]
\times\cdots\times\left[  n\right]  $).

We know that $\left(  v_{1},v_{2},\ldots,v_{n}\right)  $ is an $n$-tuple
$\left(  i_{1},i_{2},\ldots,i_{n}\right)  \in\left[  1\right]  \times\left[
2\right]  \times\cdots\times\left[  n\right]  $ such that $\sigma=t_{1,i_{1}%
}\circ t_{2,i_{2}}\circ\cdots\circ t_{n,i_{n}}$. In other words, $\left(
v_{1},v_{2},\ldots,v_{n}\right)  $ is an $n$-tuple in $\left[  1\right]
\times\left[  2\right]  \times\cdots\times\left[  n\right]  $ and satisfies
$\sigma=t_{1,v_{1}}\circ t_{2,v_{2}}\circ\cdots\circ t_{n,v_{n}}$. Thus,
$\left(  v_{1},v_{2},\ldots,v_{n}\right)  \in\left[  1\right]  \times\left[
2\right]  \times\cdots\times\left[  n\right]  $ (since $\left(  v_{1}%
,v_{2},\ldots,v_{n}\right)  $ is an $n$-tuple in $\left[  1\right]
\times\left[  2\right]  \times\cdots\times\left[  n\right]  $).

Now, $t_{1,u_{1}}\circ t_{2,u_{2}}\circ\cdots\circ t_{n,u_{n}}=\sigma
=t_{1,v_{1}}\circ t_{2,v_{2}}\circ\cdots\circ t_{n,v_{n}}$. Hence, Lemma
\ref{lem.sol.transpose.code.uni2} (applied to $k=n$) yields $\left(
u_{1},u_{2},\ldots,u_{n}\right)  =\left(  v_{1},v_{2},\ldots,v_{n}\right)  $
(since $n\in\left\{  0,1,\ldots,n\right\}  $).

Now, forget that we fixed $\left(  u_{1},u_{2},\ldots,u_{n}\right)  $ and
$\left(  v_{1},v_{2},\ldots,v_{n}\right)  $. We thus have proven that if
$\left(  u_{1},u_{2},\ldots,u_{n}\right)  $ and $\left(  v_{1},v_{2}%
,\ldots,v_{n}\right)  $ are two $n$-tuples $\left(  i_{1},i_{2},\ldots
,i_{n}\right)  \in\left[  1\right]  \times\left[  2\right]  \times\cdots
\times\left[  n\right]  $ such that $\sigma=t_{1,i_{1}}\circ t_{2,i_{2}}%
\circ\cdots\circ t_{n,i_{n}}$, then $\left(  u_{1},u_{2},\ldots,u_{n}\right)
=\left(  v_{1},v_{2},\ldots,v_{n}\right)  $. In other words, there exists
\textbf{at most one} $n$-tuple $\left(  i_{1},i_{2},\ldots,i_{n}\right)
\in\left[  1\right]  \times\left[  2\right]  \times\cdots\times\left[
n\right]  $ such that $\sigma=t_{1,i_{1}}\circ t_{2,i_{2}}\circ\cdots\circ
t_{n,i_{n}}$. This proves Claim 2.]

Combining Claim 1 with Claim 2, we conclude that there is a \textbf{unique}
$n$-tuple $\left(  i_{1},i_{2},\ldots,i_{n}\right)  \in\left[  1\right]
\times\left[  2\right]  \times\cdots\times\left[  n\right]  $ such that
$\sigma=t_{1,i_{1}}\circ t_{2,i_{2}}\circ\cdots\circ t_{n,i_{n}}$ (indeed, its
existence follows from Claim 1, whereas its uniqueness follows from Claim 2).
This solves Exercise \ref{exe.transpos.code}.
\end{proof}
\end{verlong}

\subsubsection{Some consequences}

We shall now use the result of Exercise \ref{exe.transpos.code} to derive a
few basic properties of symmetric groups. We begin with the following:

\begin{corollary}
\label{cor.transpos.code.bij}Let $n\in\mathbb{N}$. The map%
\begin{align*}
\left[  1\right]  \times\left[  2\right]  \times\cdots\times\left[  n\right]
&  \rightarrow S_{n},\\
\left(  i_{1},i_{2},\ldots,i_{n}\right)   &  \mapsto t_{1,i_{1}}\circ
t_{2,i_{2}}\circ\cdots\circ t_{n,i_{n}}%
\end{align*}
is well-defined and bijective.
\end{corollary}

\begin{vershort}
\begin{proof}
[Proof of Corollary \ref{cor.transpos.code.bij}.]For each $\left(  i_{1}%
,i_{2},\ldots,i_{n}\right)  \in\left[  1\right]  \times\left[  2\right]
\times\cdots\times\left[  n\right]  $, we have $t_{1,i_{1}}\circ t_{2,i_{2}%
}\circ\cdots\circ t_{n,i_{n}}\in S_{n}$ (since $t_{1,i_{1}},t_{2,i_{2}}%
,\ldots,t_{n,i_{n}}$ are well-defined elements of $S_{n}$). Hence, the map%
\begin{align*}
\left[  1\right]  \times\left[  2\right]  \times\cdots\times\left[  n\right]
&  \rightarrow S_{n},\\
\left(  i_{1},i_{2},\ldots,i_{n}\right)   &  \mapsto t_{1,i_{1}}\circ
t_{2,i_{2}}\circ\cdots\circ t_{n,i_{n}}%
\end{align*}
is well-defined. Let us denote this map by $A$.

The map $A$ is injective\footnote{\textit{Proof.} Let $\mathbf{x}$ and
$\mathbf{y}$ be two elements of $\left[  1\right]  \times\left[  2\right]
\times\cdots\times\left[  n\right]  $ satisfying $A\left(  \mathbf{x}\right)
=A\left(  \mathbf{y}\right)  $. We shall prove that $\mathbf{x}=\mathbf{y}$.
\par
Define $\sigma\in S_{n}$ by $\sigma=A\left(  \mathbf{x}\right)  $. Thus,
$\sigma=A\left(  \mathbf{x}\right)  =A\left(  \mathbf{y}\right)  $.
\par
Write the element $\mathbf{x}\in\left[  1\right]  \times\left[  2\right]
\times\cdots\times\left[  n\right]  $ in the form $\mathbf{x}=\left(
x_{1},x_{2},\ldots,x_{n}\right)  $. Hence, $\left(  x_{1},x_{2},\ldots
,x_{n}\right)  =\mathbf{x}\in\left[  1\right]  \times\left[  2\right]
\times\cdots\times\left[  n\right]  $. Now,%
\[
\sigma=A\left(  \underbrace{\mathbf{x}}_{=\left(  x_{1},x_{2},\ldots
,x_{n}\right)  }\right)  =A\left(  \left(  x_{1},x_{2},\ldots,x_{n}\right)
\right)  =t_{1,x_{1}}\circ t_{2,x_{2}}\circ\cdots\circ t_{n,x_{n}}%
\]
(by the definition of $A$). Thus, $\left(  x_{1},x_{2},\ldots,x_{n}\right)  $
is an $n$-tuple $\left(  i_{1},i_{2},\ldots,i_{n}\right)  \in\left[  1\right]
\times\left[  2\right]  \times\cdots\times\left[  n\right]  $ such that
$\sigma=t_{1,i_{1}}\circ t_{2,i_{2}}\circ\cdots\circ t_{n,i_{n}}$ (since
$\left(  x_{1},x_{2},\ldots,x_{n}\right)  \in\left[  1\right]  \times\left[
2\right]  \times\cdots\times\left[  n\right]  $ and $\sigma=t_{1,x_{1}}\circ
t_{2,x_{2}}\circ\cdots\circ t_{n,x_{n}}$). In other words, $\mathbf{x}$ is an
$n$-tuple $\left(  i_{1},i_{2},\ldots,i_{n}\right)  \in\left[  1\right]
\times\left[  2\right]  \times\cdots\times\left[  n\right]  $ such that
$\sigma=t_{1,i_{1}}\circ t_{2,i_{2}}\circ\cdots\circ t_{n,i_{n}}$ (since
$\mathbf{x}=\left(  x_{1},x_{2},\ldots,x_{n}\right)  $). Similarly,
$\mathbf{y}$ is an $n$-tuple $\left(  i_{1},i_{2},\ldots,i_{n}\right)
\in\left[  1\right]  \times\left[  2\right]  \times\cdots\times\left[
n\right]  $ such that $\sigma=t_{1,i_{1}}\circ t_{2,i_{2}}\circ\cdots\circ
t_{n,i_{n}}$ (since $\sigma=A\left(  \mathbf{y}\right)  $).
\par
But Exercise \ref{exe.transpos.code} shows that there is a unique $n$-tuple
$\left(  i_{1},i_{2},\ldots,i_{n}\right)  \in\left[  1\right]  \times\left[
2\right]  \times\cdots\times\left[  n\right]  $ such that $\sigma=t_{1,i_{1}%
}\circ t_{2,i_{2}}\circ\cdots\circ t_{n,i_{n}}$. In particular, there exists
\textbf{at most one} such $n$-tuple. In other words, if $\mathbf{u}$ and
$\mathbf{v}$ are two $n$-tuples $\left(  i_{1},i_{2},\ldots,i_{n}\right)
\in\left[  1\right]  \times\left[  2\right]  \times\cdots\times\left[
n\right]  $ such that $\sigma=t_{1,i_{1}}\circ t_{2,i_{2}}\circ\cdots\circ
t_{n,i_{n}}$, then $\mathbf{u}=\mathbf{v}$. Applying this to $\mathbf{u}%
=\mathbf{x}$ and $\mathbf{v}=\mathbf{y}$, we conclude that $\mathbf{x}%
=\mathbf{y}$ (since $\mathbf{x}$ and $\mathbf{y}$ are two $n$-tuples $\left(
i_{1},i_{2},\ldots,i_{n}\right)  \in\left[  1\right]  \times\left[  2\right]
\times\cdots\times\left[  n\right]  $ such that $\sigma=t_{1,i_{1}}\circ
t_{2,i_{2}}\circ\cdots\circ t_{n,i_{n}}$).
\par
Now, forget that we fixed $\mathbf{x}$ and $\mathbf{y}$. We thus have shown
that if $\mathbf{x}$ and $\mathbf{y}$ are two elements of $\left[  1\right]
\times\left[  2\right]  \times\cdots\times\left[  n\right]  $ satisfying
$A\left(  \mathbf{x}\right)  =A\left(  \mathbf{y}\right)  $, then
$\mathbf{x}=\mathbf{y}$. In other words, the map $A$ is injective.} and
surjective\footnote{\textit{Proof.} Let $\sigma\in S_{n}$. Then, Exercise
\ref{exe.transpos.code} shows that there is a unique $n$-tuple $\left(
i_{1},i_{2},\ldots,i_{n}\right)  \in\left[  1\right]  \times\left[  2\right]
\times\cdots\times\left[  n\right]  $ such that $\sigma=t_{1,i_{1}}\circ
t_{2,i_{2}}\circ\cdots\circ t_{n,i_{n}}$. Consider this $\left(  i_{1}%
,i_{2},\ldots,i_{n}\right)  $. The definition of $A$ yields $A\left(  \left(
i_{1},i_{2},\ldots,i_{n}\right)  \right)  =t_{1,i_{1}}\circ t_{2,i_{2}}%
\circ\cdots\circ t_{n,i_{n}}=\sigma$. Thus, $\sigma=A\left(
\underbrace{\left(  i_{1},i_{2},\ldots,i_{n}\right)  }_{\in\left[  1\right]
\times\left[  2\right]  \times\cdots\times\left[  n\right]  }\right)  \in
A\left(  \left[  1\right]  \times\left[  2\right]  \times\cdots\times\left[
n\right]  \right)  $.
\par
Now, forget that we fixed $\sigma$. We thus have proven that each $\sigma\in
S_{n}$ satisfies $\sigma\in A\left(  \left[  1\right]  \times\left[  2\right]
\times\cdots\times\left[  n\right]  \right)  $. In other words, $S_{n}%
\subseteq A\left(  \left[  1\right]  \times\left[  2\right]  \times
\cdots\times\left[  n\right]  \right)  $. In other words, the map $A$ is
surjective.}. Hence, the map $A$ is bijective. In other words, the map%
\begin{align*}
\left[  1\right]  \times\left[  2\right]  \times\cdots\times\left[  n\right]
&  \rightarrow S_{n},\\
\left(  i_{1},i_{2},\ldots,i_{n}\right)   &  \mapsto t_{1,i_{1}}\circ
t_{2,i_{2}}\circ\cdots\circ t_{n,i_{n}}%
\end{align*}
is bijective (since this map is precisely the map $A$). This completes the
proof of Corollary \ref{cor.transpos.code.bij}.
\end{proof}
\end{vershort}

\begin{verlong}
\begin{proof}
[Proof of Corollary \ref{cor.transpos.code.bij}.]We have $\left[  n\right]
=\left\{  1,2,\ldots,n\right\}  $ (by the definition of $\left[  n\right]  $).

For each $\left(  i_{1},i_{2},\ldots,i_{n}\right)  \in\left[  1\right]
\times\left[  2\right]  \times\cdots\times\left[  n\right]  $, we have
$t_{1,i_{1}}\circ t_{2,i_{2}}\circ\cdots\circ t_{n,i_{n}}\in S_{n}%
$\ \ \ \ \footnote{\textit{Proof.} Let $\left(  i_{1},i_{2},\ldots
,i_{n}\right)  \in\left[  1\right]  \times\left[  2\right]  \times\cdots
\times\left[  n\right]  $. We must show that $t_{1,i_{1}}\circ t_{2,i_{2}%
}\circ\cdots\circ t_{n,i_{n}}\in S_{n}$.
\par
Fix $k\in\left\{  1,2,\ldots,n\right\}  $. We have $\left(  i_{1},i_{2}%
,\ldots,i_{n}\right)  \in\left[  1\right]  \times\left[  2\right]
\times\cdots\times\left[  n\right]  $. In other words, $i_{j}\in\left[
j\right]  $ for each $j\in\left\{  1,2,\ldots,n\right\}  $. Applying this to
$j=k$, we obtain $i_{k}\in\left[  k\right]  =\left\{  1,2,\ldots,k\right\}  $
(by the definition of $\left[  k\right]  $). But $k\leq n$ (since
$k\in\left\{  1,2,\ldots,n\right\}  $) and thus $\left\{  1,2,\ldots
,k\right\}  \subseteq\left\{  1,2,\ldots,n\right\}  $. Hence, $i_{k}%
\in\left\{  1,2,\ldots,k\right\}  \subseteq\left\{  1,2,\ldots,n\right\}  $.
Hence, $k$ and $i_{k}$ are two elements of $\left\{  1,2,\ldots,n\right\}  $.
Thus, the permutation $t_{k,i_{k}}\in S_{n}$ is well-defined. In particular,
we have $t_{k,i_{k}}\in S_{n}$.
\par
Now, forget that we fixed $k$. We thus have proven that $t_{k,i_{k}}\in S_{n}$
for each $k\in\left\{  1,2,\ldots,n\right\}  $. In other words, $t_{1,i_{1}%
},t_{2,i_{2}},\ldots,t_{n,i_{n}}$ are permutations in $S_{n}$. Hence, their
composition $t_{1,i_{1}}\circ t_{2,i_{2}}\circ\cdots\circ t_{n,i_{n}}$ is a
permutation in $S_{n}$ as well. In other words, $t_{1,i_{1}}\circ t_{2,i_{2}%
}\circ\cdots\circ t_{n,i_{n}}\in S_{n}$. Qed.}. Hence, the map%
\begin{align*}
\left[  1\right]  \times\left[  2\right]  \times\cdots\times\left[  n\right]
&  \rightarrow S_{n},\\
\left(  i_{1},i_{2},\ldots,i_{n}\right)   &  \mapsto t_{1,i_{1}}\circ
t_{2,i_{2}}\circ\cdots\circ t_{n,i_{n}}%
\end{align*}
is well-defined. Let us denote this map by $A$.

The map $A$ is injective\footnote{\textit{Proof.} Let $\mathbf{x}$ and
$\mathbf{y}$ be two elements of $\left[  1\right]  \times\left[  2\right]
\times\cdots\times\left[  n\right]  $ satisfying $A\left(  \mathbf{x}\right)
=A\left(  \mathbf{y}\right)  $. We shall prove that $\mathbf{x}=\mathbf{y}$.
\par
Define $\sigma\in S_{n}$ by $\sigma=A\left(  \mathbf{x}\right)  $. Thus,
$\sigma=A\left(  \mathbf{x}\right)  =A\left(  \mathbf{y}\right)  $.
\par
Write the element $\mathbf{x}\in\left[  1\right]  \times\left[  2\right]
\times\cdots\times\left[  n\right]  $ in the form $\mathbf{x}=\left(
x_{1},x_{2},\ldots,x_{n}\right)  $. Hence, $\left(  x_{1},x_{2},\ldots
,x_{n}\right)  =\mathbf{x}\in\left[  1\right]  \times\left[  2\right]
\times\cdots\times\left[  n\right]  $.
\par
Write the element $\mathbf{y}\in\left[  1\right]  \times\left[  2\right]
\times\cdots\times\left[  n\right]  $ in the form $\mathbf{y}=\left(
y_{1},y_{2},\ldots,y_{n}\right)  $. Hence, $\left(  y_{1},y_{2},\ldots
,y_{n}\right)  =\mathbf{y}\in\left[  1\right]  \times\left[  2\right]
\times\cdots\times\left[  n\right]  $.
\par
Now,%
\[
\sigma=A\left(  \underbrace{\mathbf{x}}_{=\left(  x_{1},x_{2},\ldots
,x_{n}\right)  }\right)  =A\left(  \left(  x_{1},x_{2},\ldots,x_{n}\right)
\right)  =t_{1,x_{1}}\circ t_{2,x_{2}}\circ\cdots\circ t_{n,x_{n}}%
\]
(by the definition of $A$). Thus, $\left(  x_{1},x_{2},\ldots,x_{n}\right)  $
is an $n$-tuple $\left(  i_{1},i_{2},\ldots,i_{n}\right)  \in\left[  1\right]
\times\left[  2\right]  \times\cdots\times\left[  n\right]  $ such that
$\sigma=t_{1,i_{1}}\circ t_{2,i_{2}}\circ\cdots\circ t_{n,i_{n}}$ (since
$\left(  x_{1},x_{2},\ldots,x_{n}\right)  \in\left[  1\right]  \times\left[
2\right]  \times\cdots\times\left[  n\right]  $ and $\sigma=t_{1,x_{1}}\circ
t_{2,x_{2}}\circ\cdots\circ t_{n,x_{n}}$). In other words, $\mathbf{x}$ is an
$n$-tuple $\left(  i_{1},i_{2},\ldots,i_{n}\right)  \in\left[  1\right]
\times\left[  2\right]  \times\cdots\times\left[  n\right]  $ such that
$\sigma=t_{1,i_{1}}\circ t_{2,i_{2}}\circ\cdots\circ t_{n,i_{n}}$ (since
$\mathbf{x}=\left(  x_{1},x_{2},\ldots,x_{n}\right)  $).
\par
Also,%
\[
\sigma=A\left(  \underbrace{\mathbf{y}}_{=\left(  y_{1},y_{2},\ldots
,y_{n}\right)  }\right)  =A\left(  \left(  y_{1},y_{2},\ldots,y_{n}\right)
\right)  =t_{1,y_{1}}\circ t_{2,y_{2}}\circ\cdots\circ t_{n,y_{n}}%
\]
(by the definition of $A$). Thus, $\left(  y_{1},y_{2},\ldots,y_{n}\right)  $
is an $n$-tuple $\left(  i_{1},i_{2},\ldots,i_{n}\right)  \in\left[  1\right]
\times\left[  2\right]  \times\cdots\times\left[  n\right]  $ such that
$\sigma=t_{1,i_{1}}\circ t_{2,i_{2}}\circ\cdots\circ t_{n,i_{n}}$ (since
$\left(  y_{1},y_{2},\ldots,y_{n}\right)  \in\left[  1\right]  \times\left[
2\right]  \times\cdots\times\left[  n\right]  $ and $\sigma=t_{1,y_{1}}\circ
t_{2,y_{2}}\circ\cdots\circ t_{n,y_{n}}$). In other words, $\mathbf{y}$ is an
$n$-tuple $\left(  i_{1},i_{2},\ldots,i_{n}\right)  \in\left[  1\right]
\times\left[  2\right]  \times\cdots\times\left[  n\right]  $ such that
$\sigma=t_{1,i_{1}}\circ t_{2,i_{2}}\circ\cdots\circ t_{n,i_{n}}$ (since
$\mathbf{y}=\left(  y_{1},y_{2},\ldots,y_{n}\right)  $).
\par
But Exercise \ref{exe.transpos.code} shows that there is a unique $n$-tuple
$\left(  i_{1},i_{2},\ldots,i_{n}\right)  \in\left[  1\right]  \times\left[
2\right]  \times\cdots\times\left[  n\right]  $ such that $\sigma=t_{1,i_{1}%
}\circ t_{2,i_{2}}\circ\cdots\circ t_{n,i_{n}}$. In particular, there exists
\textbf{at most one} such $n$-tuple. In other words, if $\mathbf{u}$ and
$\mathbf{v}$ are two $n$-tuples $\left(  i_{1},i_{2},\ldots,i_{n}\right)
\in\left[  1\right]  \times\left[  2\right]  \times\cdots\times\left[
n\right]  $ such that $\sigma=t_{1,i_{1}}\circ t_{2,i_{2}}\circ\cdots\circ
t_{n,i_{n}}$, then $\mathbf{u}=\mathbf{v}$. Applying this to $\mathbf{u}%
=\mathbf{x}$ and $\mathbf{v}=\mathbf{y}$, we conclude that $\mathbf{x}%
=\mathbf{y}$ (since $\mathbf{x}$ and $\mathbf{y}$ are two $n$-tuples $\left(
i_{1},i_{2},\ldots,i_{n}\right)  \in\left[  1\right]  \times\left[  2\right]
\times\cdots\times\left[  n\right]  $ such that $\sigma=t_{1,i_{1}}\circ
t_{2,i_{2}}\circ\cdots\circ t_{n,i_{n}}$).
\par
Now, forget that we fixed $\mathbf{x}$ and $\mathbf{y}$. We thus have shown
that if $\mathbf{x}$ and $\mathbf{y}$ are two elements of $\left[  1\right]
\times\left[  2\right]  \times\cdots\times\left[  n\right]  $ satisfying
$A\left(  \mathbf{x}\right)  =A\left(  \mathbf{y}\right)  $, then
$\mathbf{x}=\mathbf{y}$. In other words, the map $A$ is injective.} and
surjective\footnote{\textit{Proof.} Let $\sigma\in S_{n}$. Then, Exercise
\ref{exe.transpos.code} shows that there is a unique $n$-tuple $\left(
i_{1},i_{2},\ldots,i_{n}\right)  \in\left[  1\right]  \times\left[  2\right]
\times\cdots\times\left[  n\right]  $ such that $\sigma=t_{1,i_{1}}\circ
t_{2,i_{2}}\circ\cdots\circ t_{n,i_{n}}$. Consider this $\left(  i_{1}%
,i_{2},\ldots,i_{n}\right)  $. The definition of $A$ yields $A\left(  \left(
i_{1},i_{2},\ldots,i_{n}\right)  \right)  =t_{1,i_{1}}\circ t_{2,i_{2}}%
\circ\cdots\circ t_{n,i_{n}}=\sigma$. Thus, $\sigma=A\left(
\underbrace{\left(  i_{1},i_{2},\ldots,i_{n}\right)  }_{\in\left[  1\right]
\times\left[  2\right]  \times\cdots\times\left[  n\right]  }\right)  \in
A\left(  \left[  1\right]  \times\left[  2\right]  \times\cdots\times\left[
n\right]  \right)  $.
\par
Now, forget that we fixed $\sigma$. We thus have proven that each $\sigma\in
S_{n}$ satisfies $\sigma\in A\left(  \left[  1\right]  \times\left[  2\right]
\times\cdots\times\left[  n\right]  \right)  $. In other words, $S_{n}%
\subseteq A\left(  \left[  1\right]  \times\left[  2\right]  \times
\cdots\times\left[  n\right]  \right)  $. In other words, the map $A$ is
surjective.}. Hence, the map $A$ is bijective. In other words, the map%
\begin{align*}
\left[  1\right]  \times\left[  2\right]  \times\cdots\times\left[  n\right]
&  \rightarrow S_{n},\\
\left(  i_{1},i_{2},\ldots,i_{n}\right)   &  \mapsto t_{1,i_{1}}\circ
t_{2,i_{2}}\circ\cdots\circ t_{n,i_{n}}%
\end{align*}
is bijective\footnote{since $A$ is the map%
\begin{align*}
\left[  1\right]  \times\left[  2\right]  \times\cdots\times\left[  n\right]
&  \rightarrow S_{n},\\
\left(  i_{1},i_{2},\ldots,i_{n}\right)   &  \mapsto t_{1,i_{1}}\circ
t_{2,i_{2}}\circ\cdots\circ t_{n,i_{n}}%
\end{align*}
(by the definition of $A$)}. This completes the proof of Corollary
\ref{cor.transpos.code.bij}.
\end{proof}
\end{verlong}

We can use Corollary \ref{cor.transpos.code.bij} to obtain the following:

\begin{corollary}
\label{cor.transpos.code.n!}Let $n\in\mathbb{N}$. Then, $\left\vert
S_{n}\right\vert =n!$.
\end{corollary}

\begin{proof}
[Proof of Corollary \ref{cor.transpos.code.n!}.]Corollary
\ref{cor.transpos.code.bij} shows that the map%
\begin{align*}
\left[  1\right]  \times\left[  2\right]  \times\cdots\times\left[  n\right]
&  \rightarrow S_{n},\\
\left(  i_{1},i_{2},\ldots,i_{n}\right)   &  \mapsto t_{1,i_{1}}\circ
t_{2,i_{2}}\circ\cdots\circ t_{n,i_{n}}%
\end{align*}
is well-defined and bijective. Hence, this map is a bijection from $\left[
1\right]  \times\left[  2\right]  \times\cdots\times\left[  n\right]  $ to
$S_{n}$. Thus, there exists a bijection from $\left[  1\right]  \times\left[
2\right]  \times\cdots\times\left[  n\right]  $ to $S_{n}$ (namely, this map).
Thus,%
\begin{align*}
\left\vert S_{n}\right\vert  &  =\left\vert \left[  1\right]  \times\left[
2\right]  \times\cdots\times\left[  n\right]  \right\vert =\left\vert \left[
1\right]  \right\vert \cdot\left\vert \left[  2\right]  \right\vert
\cdot\cdots\cdot\left\vert \left[  n\right]  \right\vert =\prod_{k=1}%
^{n}\left\vert \underbrace{\left[  k\right]  }_{\substack{=\left\{
1,2,\ldots,k\right\}  \\\text{(by the definition of }\left[  k\right]
\text{)}}}\right\vert \\
&  =\prod_{k=1}^{n}\underbrace{\left\vert \left\{  1,2,\ldots,k\right\}
\right\vert }_{=k}=\prod_{k=1}^{n}k=1\cdot2\cdot\cdots\cdot n=n!.
\end{align*}
This proves Corollary \ref{cor.transpos.code.n!}.
\end{proof}

Corollary \ref{cor.transpos.code.n!} can be generalized:

\begin{corollary}
\label{cor.transpose.code.X!}Let $X$ be a finite set. Then, the number of all
permutations of $X$ is $\left\vert X\right\vert !$.
\end{corollary}

To derive Corollary \ref{cor.transpose.code.X!} from Corollary
\ref{cor.transpos.code.n!}, we need a basic lemma:

\begin{lemma}
\label{lem.transpose.code.XY}Let $X$ and $Y$ be two sets. Let $f:X\rightarrow
Y$ be a bijection. Then, the map%
\begin{align*}
\left\{  \text{permutations of }X\right\}   &  \rightarrow\left\{
\text{permutations of }Y\right\}  ,\\
\sigma &  \mapsto f\circ\sigma\circ f^{-1}%
\end{align*}
is well-defined and bijective.
\end{lemma}

\begin{vershort}
\begin{proof}
[Proof of Lemma \ref{lem.transpose.code.XY}.]This is straightforward to check.
(The inverse of this map is the map%
\begin{align*}
\left\{  \text{permutations of }Y\right\}   &  \rightarrow\left\{
\text{permutations of }X\right\}  ,\\
\tau &  \mapsto f^{-1}\circ\tau\circ f.
\end{align*}
)
\end{proof}
\end{vershort}

\begin{verlong}
\begin{proof}
[Proof of Lemma \ref{lem.transpose.code.XY}.]The map $f:X\rightarrow Y$ is a
bijection. Thus, its inverse $f^{-1}:Y\rightarrow X$ is well-defined and is a
bijection as well.

For each $\sigma\in\left\{  \text{permutations of }X\right\}  $, we have
$f\circ\sigma\circ f^{-1}\in\left\{  \text{permutations of }Y\right\}
$\ \ \ \ \footnote{\textit{Proof.} Let $\sigma\in\left\{  \text{permutations
of }X\right\}  $. We must prove that $f\circ\sigma\circ f^{-1}\in\left\{
\text{permutations of }Y\right\}  $.
\par
We have $\sigma\in\left\{  \text{permutations of }X\right\}  $. In other
words, $\sigma$ is a permutation of $X$. In other words, $\sigma$ is a
bijection $X\rightarrow X$.
\par
The three maps $f$, $\sigma$ and $f^{-1}$ are bijections. Thus, their
composition $f\circ\sigma\circ f^{-1}$ is a bijection as well. Hence,
$f\circ\sigma\circ f^{-1}$ is a bijection $Y\rightarrow Y$. In other words,
$f\circ\sigma\circ f^{-1}$ is a permutation of $Y$. In other words,
$f\circ\sigma\circ f^{-1}\in\left\{  \text{permutations of }Y\right\}  $.
Qed.}. Hence, the map%
\begin{align*}
\left\{  \text{permutations of }X\right\}   &  \rightarrow\left\{
\text{permutations of }Y\right\}  ,\\
\sigma &  \mapsto f\circ\sigma\circ f^{-1}%
\end{align*}
is well-defined. Let us denote this map by $A$.

For each $\tau\in\left\{  \text{permutations of }Y\right\}  $, we have
$f^{-1}\circ\tau\circ f\in\left\{  \text{permutations of }X\right\}
$\ \ \ \ \footnote{\textit{Proof.} Let $\tau\in\left\{  \text{permutations of
}Y\right\}  $. We must prove that $f^{-1}\circ\tau\circ f\in\left\{
\text{permutations of }X\right\}  $.
\par
We have $\tau\in\left\{  \text{permutations of }Y\right\}  $. In other words,
$\tau$ is a permutation of $Y$. In other words, $\tau$ is a bijection
$Y\rightarrow Y$.
\par
The three maps $f^{-1}$, $\tau$ and $f$ are bijections. Thus, their
composition $f^{-1}\circ\tau\circ f$ is a bijection as well. Hence,
$f^{-1}\circ\tau\circ f$ is a bijection $X\rightarrow X$. In other words,
$f^{-1}\circ\tau\circ f$ is a permutation of $X$. In other words, $f^{-1}%
\circ\tau\circ f\in\left\{  \text{permutations of }X\right\}  $. Qed.}. Hence,
the map%
\begin{align*}
\left\{  \text{permutations of }Y\right\}   &  \rightarrow\left\{
\text{permutations of }X\right\}  ,\\
\tau &  \mapsto f^{-1}\circ\tau\circ f
\end{align*}
is well-defined. Let us denote this map by $B$.

We have $A\circ B=\operatorname*{id}$\ \ \ \ \footnote{\textit{Proof.} Let
$\tau\in\left\{  \text{permutations of }Y\right\}  $. Then, the definition of
$B$ yields $B\left(  \tau\right)  =f^{-1}\circ\tau\circ f$. But the definition
of $A$ yields%
\[
A\left(  B\left(  \tau\right)  \right)  =f\circ\underbrace{\left(  B\left(
\tau\right)  \right)  }_{=f^{-1}\circ\tau\circ f}\circ f^{-1}%
=\underbrace{f\circ f^{-1}}_{=\operatorname*{id}}\circ\tau\circ
\underbrace{f\circ f^{-1}}_{=\operatorname*{id}}=\tau.
\]
Hence, $\left(  A\circ B\right)  \left(  \tau\right)  =A\left(  B\left(
\tau\right)  \right)  =\tau=\operatorname*{id}\left(  \tau\right)  $.
\par
Now, forget that we fixed $\tau$. We thus have shown that $\left(  A\circ
B\right)  \left(  \tau\right)  =\operatorname*{id}\left(  \tau\right)  $ for
each $\tau\in\left\{  \text{permutations of }Y\right\}  $. In other words,
$A\circ B=\operatorname*{id}$.} and $B\circ A=\operatorname*{id}%
$\ \ \ \ \footnote{\textit{Proof.} Let $\sigma\in\left\{  \text{permutations
of }X\right\}  $. Then, the definition of $A$ yields $A\left(  \sigma\right)
=f\circ\sigma\circ f^{-1}$. But the definition of $B$ yields%
\[
B\left(  A\left(  \sigma\right)  \right)  =f^{-1}\circ\underbrace{\left(
A\left(  \sigma\right)  \right)  }_{=f\circ\sigma\circ f^{-1}}\circ
f=\underbrace{f^{-1}\circ f}_{=\operatorname*{id}}\circ\sigma\circ
\underbrace{f^{-1}\circ f}_{=\operatorname*{id}}=\sigma.
\]
Hence, $\left(  B\circ A\right)  \left(  \sigma\right)  =B\left(  A\left(
\sigma\right)  \right)  =\sigma=\operatorname*{id}\left(  \sigma\right)  $.
\par
Now, forget that we fixed $\sigma$. We thus have shown that $\left(  B\circ
A\right)  \left(  \sigma\right)  =\operatorname*{id}\left(  \sigma\right)  $
for each $\sigma\in\left\{  \text{permutations of }X\right\}  $. In other
words, $B\circ A=\operatorname*{id}$.}. These two equalities show that the
maps $A$ and $B$ are mutually inverse. Thus, the map $A$ is invertible. In
other words, the map $A$ is bijective. In other words, the map%
\begin{align*}
\left\{  \text{permutations of }X\right\}   &  \rightarrow\left\{
\text{permutations of }Y\right\}  ,\\
\sigma &  \mapsto f\circ\sigma\circ f^{-1}%
\end{align*}
is bijective\footnote{since the map $A$ is the map
\begin{align*}
\left\{  \text{permutations of }X\right\}   &  \rightarrow\left\{
\text{permutations of }Y\right\}  ,\\
\sigma &  \mapsto f\circ\sigma\circ f^{-1}%
\end{align*}
(by the definition of $A$)}. This completes the proof of Lemma
\ref{lem.transpose.code.XY}.
\end{proof}
\end{verlong}

\begin{proof}
[Proof of Corollary \ref{cor.transpose.code.X!}.]Define $n\in\mathbb{N}$ by
$n=\left\vert X\right\vert $. (This is well-defined, since $X$ is a finite set.)

The definition of the symmetric group $S_{n}$ shows that $S_{n}$ is the set of
all permutations of the set $\left\{  1,2,\ldots,n\right\}  $. In other words,%
\begin{equation}
S_{n}=\left\{  \text{permutations of the set }\left\{  1,2,\ldots,n\right\}
\right\}  . \label{pf.cor.transpose.code.X!.1}%
\end{equation}

Define a finite set $Y$ by $Y=\left\{  1,2,\ldots,n\right\}  $. Then,
$\left\vert \underbrace{Y}_{=\left\{  1,2,\ldots,n\right\}  }\right\vert
=\left\vert \left\{  1,2,\ldots,n\right\}  \right\vert =n$. Also,%
\begin{align}
&  \left\{  \text{permutations of }\underbrace{Y}_{=\left\{  1,2,\ldots
,n\right\}  }\right\} \nonumber\\
&  =\left\{  \text{permutations of }\left\{  1,2,\ldots,n\right\}  \right\}
=\left\{  \text{permutations of the set }\left\{  1,2,\ldots,n\right\}
\right\} \nonumber\\
&  =S_{n}\ \ \ \ \ \ \ \ \ \ \left(  \text{by
(\ref{pf.cor.transpose.code.X!.1})}\right)  .
\label{pf.cor.transpose.code.X!.2}%
\end{align}

But we have $\left\vert Y\right\vert =n=\left\vert X\right\vert $. In other
words, the two finite sets $X$ and $Y$ have the same size. Hence, there exists
a bijection $f:X\rightarrow Y$. Consider this $f$. Lemma
\ref{lem.transpose.code.XY} thus shows that the map%
\begin{align*}
\left\{  \text{permutations of }X\right\}   &  \rightarrow\left\{
\text{permutations of }Y\right\}  ,\\
\sigma &  \mapsto f\circ\sigma\circ f^{-1}%
\end{align*}
is well-defined and bijective. Thus, this map is a bijection. Hence, there
exists a bijection $\left\{  \text{permutations of }X\right\}  \rightarrow
\left\{  \text{permutations of }Y\right\}  $ (namely, this map). Thus,%
\[
\left\vert \left\{  \text{permutations of }X\right\}  \right\vert =\left\vert
\underbrace{\left\{  \text{permutations of }Y\right\}  }_{\substack{=S_{n}%
\\\text{(by (\ref{pf.cor.transpose.code.X!.2}))}}}\right\vert =\left\vert
S_{n}\right\vert =n!
\]
(by Corollary \ref{cor.transpos.code.n!}). Now, the number of all permutations
of $X$ is $\left\vert \left\{  \text{permutations of }X\right\}  \right\vert
=\underbrace{n}_{=\left\vert X\right\vert }!=\left\vert X\right\vert !$. This
proves Corollary \ref{cor.transpose.code.X!}.
\end{proof}

\subsection{Solution to Exercise \ref{exe.ps4.1ab}}

\begin{proof}
[Solution to Exercise \ref{exe.ps4.1ab}.]\textbf{(a)} We shall prove that%
\begin{equation}
\ell\left(  t_{i,j}\right)  =2\left\vert j-i\right\vert -1.
\label{sol.ps4.1ab.a.claim}%
\end{equation}

[\textit{Proof of (\ref{sol.ps4.1ab.a.claim}):} We know that $t_{i,j}$ is the
permutation in $S_{n}$ which swaps $i$ with $j$ while leaving all other
elements of $\left\{  1,2,\ldots,n\right\}  $ unchanged; on the other hand,
$t_{j,i}$ is the permutation in $S_{n}$ which swaps $j$ with $i$ while leaving
all other elements of $\left\{  1,2,\ldots,n\right\}  $ unchanged. Comparing
these two descriptions of $t_{i,j}$ and $t_{j,i}$, we immediately see that
they are identical (since swapping $i$ with $j$ is the same thing as swapping
$j$ with $i$). Thus, $t_{i,j}=t_{j,i}$. Also, clearly, $\left\vert
j-i\right\vert =\left\vert i-j\right\vert $. Hence, the claim
(\ref{sol.ps4.1ab.a.claim}) does not change if we swap $i$ with $j$. Thus, we
can WLOG assume that $i\leq j$ (because otherwise, we can just swap $i$ with
$j$). Assume this. Now, $i\leq j$, so that $i<j$ (since $i$ and $j$ are
distinct). Hence, $j>i$, so that $j-i>0$, so that $\left\vert j-i\right\vert
=j-i$.

\begin{vershort}
Set%
\begin{align*}
A  &  =\left\{  \left(  i,k\right)  \ \mid\ k\in\left\{  i+1,i+2,\ldots
,j\right\}  \right\}  \ \ \ \ \ \ \ \ \ \ \text{and}\\
B  &  =\left\{  \left(  k,j\right)  \ \mid\ k\in\left\{  i+1,i+2,\ldots
,j-1\right\}  \right\}  .
\end{align*}
These two sets $A$ and $B$ satisfy $\left\vert A\right\vert =j-i$ and
$\left\vert B\right\vert =j-i-1$. Also, it is easy to see that these sets $A$
and $B$ are disjoint\footnote{\textit{Proof.} Every element of $B$ is a pair
$\left(  k,j\right)  $ whose first entry is $>i$, whereas every element of $A$
is a pair $\left(  i,k\right)  $ whose first entry equals $i$. Thus, if the
sets $A$ and $B$ had an element $e$ in common, then $e$ would be a pair whose
first entry is $>i$ (since $e\in B$) and equals $i$ (since $e\in A$) at the
same time, which of course is impossible. Hence, $A$ and $B$ are disjoint.}.
Thus,%
\[
\left\vert A\cup B\right\vert =\underbrace{\left\vert A\right\vert }%
_{=j-i}+\underbrace{\left\vert B\right\vert }_{=j-i-1}=\left(  j-i\right)
+\left(  j-i-1\right)  =2\underbrace{\left(  j-i\right)  }_{=\left\vert
j-i\right\vert }-1=2\left\vert j-i\right\vert -1.
\]

Now, let $\operatorname*{Inv}\left(  t_{i,j}\right)  $ denote the set of all
inversions of $t_{i,j}$. Then, $\ell\left(  t_{i,j}\right)  =\left\vert
\operatorname*{Inv}\left(  t_{i,j}\right)  \right\vert $ (because $\ell\left(
t_{i,j}\right)  $ was defined as the number of inversions of $t_{i,j}$, which
number is obviously $\left\vert \operatorname*{Inv}\left(  t_{i,j}\right)
\right\vert $).

We shall now show that $\operatorname*{Inv}\left(  t_{i,j}\right)  =A\cup B$.
Indeed, it is clearly enough to prove $A\cup B\subseteq\operatorname*{Inv}%
\left(  t_{i,j}\right)  $ and $\operatorname*{Inv}\left(  t_{i,j}\right)
\subseteq A\cup B$. Proving that $A\cup B\subseteq\operatorname*{Inv}\left(
t_{i,j}\right)  $ means proving that every element of $A\cup B$ is an
inversion of $t_{i,j}$; this is straightforward\footnote{\textit{Proof.} We
want to show that $A\cup B\subseteq\operatorname*{Inv}\left(  t_{i,j}\right)
$. In other words, we want to prove that $e\in\operatorname*{Inv}\left(
t_{i,j}\right)  $ for every $e\in A\cup B$.
\par
So let $e\in A\cup B$. Thus, either $e\in A$ or $e\in B$.
\par
Let us first consider the case when $e\in A$. Thus, $e\in A=\left\{  \left(
i,k\right)  \ \mid\ k\in\left\{  i+1,i+2,\ldots,j\right\}  \right\}  $. In
other words, $e$ has the form $e=\left(  i,k\right)  $ for some $k\in\left\{
i+1,i+2,\ldots,j\right\}  $. Consider this $k$. The permutation $t_{i,j}$
swaps $i$ with $j$ while leaving all other numbers fixed. Thus, the
permutation $t_{i,j}$ leaves the numbers $i+1,i+2,\ldots,j-1$ fixed, while
sending the number $j$ to $i$. Consequently, $t_{i,j}$ sends the numbers
$i+1,i+2,\ldots,j-1,j$ to $i+1,i+2,\ldots,j-1,i$, respectively. Notice that
all of the latter numbers $i+1,i+2,\ldots,j-1,i$ are smaller than $j$. Thus,
$t_{i,j}\left(  p\right)  <j$ for every $p\in\left\{  i+1,i+2,\ldots
,j\right\}  $. Applying this to $p=k$, we conclude that $t_{i,j}\left(
k\right)  <j$.
\par
We have $i<k$ (since $k\in\left\{  i+1,i+2,\ldots,j\right\}  $) but
$t_{i,j}\left(  i\right)  =j>t_{i,j}\left(  k\right)  $ (since we have just
showed that $t_{i,j}\left(  k\right)  <j$). Thus, $\left(  i,k\right)  $ is an
inversion of $t_{i,j}$. In other words, $\left(  i,k\right)  \in
\operatorname*{Inv}\left(  t_{i,j}\right)  $. Thus, $e=\left(  i,k\right)
\in\operatorname*{Inv}\left(  t_{i,j}\right)  $.
\par
Thus, we have proven that $e\in\operatorname*{Inv}\left(  t_{i,j}\right)  $ in
the case when $e\in A$. A similar argument (but now using $t_{i,j}\left(
k\right)  >i$ instead of $t_{i,j}\left(  k\right)  <j$) shows that
$e\in\operatorname*{Inv}\left(  t_{i,j}\right)  $ in the case when $e\in B$.
Since either of these two cases must hold (because we have either $e\in A$ or
$e\in B$), we thus conclude that $e\in\operatorname*{Inv}\left(
t_{i,j}\right)  $. This concludes the proof.}. Proving that
$\operatorname*{Inv}\left(  t_{i,j}\right)  \subseteq A\cup B$ means proving
that every inversion of $t_{i,j}$ belongs to $A\cup B$; this is equally
straightforward (although more tiresome)\footnote{\textit{Proof.} We want to
show that $\operatorname*{Inv}\left(  t_{i,j}\right)  \subseteq A\cup B$. In
other words, we want to prove that $c\in A\cup B$ for every $c\in
\operatorname*{Inv}\left(  t_{i,j}\right)  $.
\par
So let $c\in\operatorname*{Inv}\left(  t_{i,j}\right)  $. Thus, $c$ is an
inversion of $t_{i,j}$. In other words, $c$ is a pair $\left(  u,v\right)  $
of integers satisfying $1\leq u<v\leq n$ and $t_{i,j}\left(  u\right)
>t_{i,j}\left(  v\right)  $. Consider this $\left(  u,v\right)  $. Thus,
$c=\left(  u,v\right)  $. Our goal is to show that $c\in A\cup B$.
\par
The permutation $t_{i,j}$ swaps $i$ with $j$ while leaving all other numbers
fixed. It thus makes sense to analyze several cases separately, depending on
which of the numbers $u$ and $v$ belongs to $\left\{  i,j\right\}  $. Four
cases are possible:
\par
\textit{Case 1:} We have $u\in\left\{  i,j\right\}  $ and $v\in\left\{
i,j\right\}  $.
\par
\textit{Case 2:} We have $u\in\left\{  i,j\right\}  $ and $v\notin\left\{
i,j\right\}  $.
\par
\textit{Case 3:} We have $u\notin\left\{  i,j\right\}  $ and $v\in\left\{
i,j\right\}  $.
\par
\textit{Case 4:} We have $u\notin\left\{  i,j\right\}  $ and $v\notin\left\{
i,j\right\}  $.
\par
Let us first consider Case 1. In this case, we have $u\in\left\{  i,j\right\}
$ and $v\in\left\{  i,j\right\}  $. Thus, $u$ and $v$ are two elements of
$\left\{  i,j\right\}  $. Since $u<v$, this leaves only one possibility for
the pair $\left(  u,v\right)  $: namely, $\left(  u,v\right)  =\left(
i,j\right)  $. Thus,
\begin{align*}
c  &  =\left(  u,v\right)  =\left(  i,j\right)  \in\left\{  \left(
i,k\right)  \ \mid\ k\in\left\{  i+1,i+2,\ldots,j\right\}  \right\}
\ \ \ \ \ \ \ \ \ \ \left(  \text{since }j\in\left\{  i+1,i+2,\ldots
,j\right\}  \right) \\
&  =A\subseteq A\cup B.
\end{align*}
Thus, $c\in A\cup B$ is proven in Case 1.
\par
Let us next consider Case 2. In this case, we have $u\in\left\{  i,j\right\}
$ and $v\notin\left\{  i,j\right\}  $. Since $v\notin\left\{  i,j\right\}  $,
we have $t_{i,j}\left(  v\right)  =v$ (since $t_{i,j}$ leaves all numbers
other than $i$ and $j$ unchanged). Thus, $t_{i,j}\left(  u\right)
>t_{i,j}\left(  v\right)  =v>u$ (since $u<v$). If we had $u=j$, then this
would rewrite as $t_{i,j}\left(  j\right)  >j$, which would contradict
$t_{i,j}\left(  j\right)  =i<j$. Thus, we cannot have $u=j$. Hence, we must
have $u=i$ (since $u\in\left\{  i,j\right\}  $ forces $u$ to be either $i$ or
$j$). But we have shown that $t_{i,j}\left(  u\right)  >v$, so that
$v<t_{i,j}\left(  \underbrace{u}_{=i}\right)  =t_{i,j}\left(  i\right)  =j$
(since $t_{i,j}$ swaps $i$ with $j$). Combined with $v>u=i$, this yields
$i<v<j$, so that $v\in\left\{  i+1,i+2,\ldots,j-1\right\}  $, so that%
\begin{align*}
c  &  =\left(  \underbrace{u}_{=i},v\right)  =\left(  i,v\right)  \in\left\{
\left(  i,k\right)  \ \mid\ k\in\left\{  i+1,i+2,\ldots,j-1\right\}  \right\}
\\
&  \ \ \ \ \ \ \ \ \ \ \left(  \text{since }v\in\left\{  i+1,i+2,\ldots
,j-1\right\}  \right) \\
&  \subseteq\left\{  \left(  i,k\right)  \ \mid\ k\in\left\{  i+1,i+2,\ldots
,j\right\}  \right\}  =A\subseteq A\cup B.
\end{align*}
Thus, $c\in A\cup B$ is proven in Case 2.
\par
Let us next consider Case 3. In this case, we have $u\notin\left\{
i,j\right\}  $ and $v\in\left\{  i,j\right\}  $. Since $u\notin\left\{
i,j\right\}  $, we have $t_{i,j}\left(  u\right)  =u$ (since $t_{i,j}$ leaves
all numbers other than $i$ and $j$ unchanged). Thus, from $t_{i,j}\left(
u\right)  >t_{i,j}\left(  v\right)  $, we obtain $t_{i,j}\left(  v\right)
<t_{i,j}\left(  u\right)  =u<v$. If we had $v=i$, then this would rewrite as
$t_{i,j}\left(  i\right)  <i$, which would contradict $t_{i,j}\left(
i\right)  =j>i$. Thus, we cannot have $v=i$. Hence, we must have $v=j$ (since
$v\in\left\{  i,j\right\}  $ forces $v$ to be either $i$ or $j$). But we have
shown that $t_{i,j}\left(  v\right)  <u$, so that $u>t_{i,j}\left(
\underbrace{v}_{=j}\right)  =t_{i,j}\left(  j\right)  =i$ (since $t_{i,j}$
swaps $i$ with $j$). Combined with $u<v=j$, this yields $i<u<j$, so that
$u\in\left\{  i+1,i+2,\ldots,j-1\right\}  $, so that%
\begin{align*}
c  &  =\left(  u,\underbrace{v}_{=j}\right)  =\left(  u,j\right)  \in\left\{
\left(  k,j\right)  \ \mid\ k\in\left\{  i+1,i+2,\ldots,j-1\right\}  \right\}
\\
&  \ \ \ \ \ \ \ \ \ \ \left(  \text{since }u\in\left\{  i+1,i+2,\ldots
,j-1\right\}  \right) \\
&  =B\subseteq A\cup B.
\end{align*}
Thus, $c\in A\cup B$ is proven in Case 3.
\par
(Notice that Case 3 was very similar to Case 2 -- almost like a mirror version
of that case, if not for a slight asymmetry in our definition of the sets $A$
and $B$.)
\par
Let us finally consider Case 4. In this case, we have $u\notin\left\{
i,j\right\}  $ and $v\notin\left\{  i,j\right\}  $. Thus, $t_{i,j}\left(
u\right)  =u$ (as in Case 3) and $t_{i,j}\left(  v\right)  =v$ (as in Case 2),
so that $t_{i,j}\left(  u\right)  =u<v=t_{i,j}\left(  v\right)  $. This
contradicts $t_{i,j}\left(  u\right)  >t_{i,j}\left(  v\right)  $. This
contradiction shows that Case 4 cannot happen; thus, we can ignore this case
completely. (Or we can argue that because \textquotedblleft ex falso
quodlibet\textquotedblright, we have $c\in A\cup B$ in Case 4.
\par
[The principle of \textquotedblleft ex falso quodlibet\textquotedblright%
\ (which we have already stated in Convention \ref{conv.logic.vacuous}) says
that from a false assertion, any arbitrary assertion follows. (For example, if
$1=0$, then anything is true.) This is one of the basic principles in logic,
and we could use it here to prove $c\in A\cup B$ in Case 4: Namely, since we
have derived a contradiction (i.e., proven a false assertion) in Case 4, we
see that any arbitrary assertion holds in Case 4; in particular, $c\in A\cup
B$ holds in Case 4.])
\par
We have now checked that $c\in A\cup B$ in each of the four cases 1, 2, 3 and
4 (or in each of the three cases 1, 2 and 3, and Case 4 never happens). Thus,
$c\in A\cup B$ always holds. This completes our proof.}. Hence, both $A\cup
B\subseteq\operatorname*{Inv}\left(  t_{i,j}\right)  $ and
$\operatorname*{Inv}\left(  t_{i,j}\right)  \subseteq A\cup B$ are proven, and
we conclude that $\operatorname*{Inv}\left(  t_{i,j}\right)  =A\cup B$. Thus,
$\left\vert \operatorname*{Inv}\left(  t_{i,j}\right)  \right\vert =\left\vert
A\cup B\right\vert =2\left\vert j-i\right\vert -1$. Thus, $\ell\left(
t_{i,j}\right)  =\left\vert \operatorname*{Inv}\left(  t_{i,j}\right)
\right\vert =2\left\vert j-i\right\vert -1$. This proves
(\ref{sol.ps4.1ab.a.claim}).]
\end{vershort}

\begin{verlong}
The permutation $t_{i,j}$ swaps $i$ with $j$ while leaving all other elements
of $\left\{  1,2,\ldots,n\right\}  $ unchanged. In other words, we have
$t_{i,j}\left(  i\right)  =j$, $t_{i,j}\left(  j\right)  =i$ and%
\begin{equation}
t_{i,j}\left(  k\right)  =k\ \ \ \ \ \ \ \ \ \ \text{for every }k\in\left\{
1,2,\ldots,n\right\}  \text{ satisfying }k\notin\left\{  i,j\right\}  .
\label{sol.ps4.1ab.a.tij}%
\end{equation}

Let%
\begin{align*}
A  &  =\left\{  \left(  i,k\right)  \ \mid\ k\in\left\{  i+1,i+2,\ldots
,j\right\}  \right\}  \ \ \ \ \ \ \ \ \ \ \text{and}\\
B  &  =\left\{  \left(  k,j\right)  \ \mid\ k\in\left\{  i+1,i+2,\ldots
,j-1\right\}  \right\}  .
\end{align*}
These two sets $A$ and $B$ satisfy $\left\vert A\right\vert =j-i$%
\ \ \ \ \footnote{\textit{Proof.} The elements $\left(  i,k\right)  $ for all
$k\in\left\{  i+1,i+2,\ldots,j\right\}  $ are pairwise distinct (because $k$
can be reconstructed from $\left(  i,k\right)  $). Therefore, the number of
these elements (counted without multiplicities) is $j-i$ (since the number of
all $k\in\left\{  i+1,i+2,\ldots,j\right\}  $ is $j-i$). In other words,
$\left\vert \left\{  \left(  i,k\right)  \ \mid\ k\in\left\{  i+1,i+2,\ldots
,j\right\}  \right\}  \right\vert =j-i$. Now,%
\[
\left\vert \underbrace{A}_{=\left\{  \left(  i,k\right)  \ \mid\ k\in\left\{
i+1,i+2,\ldots,j\right\}  \right\}  }\right\vert =\left\vert \left\{  \left(
i,k\right)  \ \mid\ k\in\left\{  i+1,i+2,\ldots,j\right\}  \right\}
\right\vert =j-i,
\]
qed.} and $\left\vert B\right\vert =j-i-1$\ \ \ \ \footnote{\textit{Proof.}
The elements $\left(  k,j\right)  $ for all $k\in\left\{  i+1,i+2,\ldots
,j-1\right\}  $ are pairwise distinct (because $k$ can be reconstructed from
$\left(  k,j\right)  $). Therefore, the number of these elements (counted
without multiplicities) is $j-i-1$ (since the number of all $k\in\left\{
i+1,i+2,\ldots,j-1\right\}  $ is $\left(  j-1\right)  -i=j-i-1$). In other
words, $\left\vert \left\{  \left(  k,j\right)  \ \mid\ k\in\left\{
i+1,i+2,\ldots,j-1\right\}  \right\}  \right\vert =j-i-1$. Now,%
\[
\left\vert \underbrace{B}_{=\left\{  \left(  k,j\right)  \ \mid\ k\in\left\{
i+1,i+2,\ldots,j-1\right\}  \right\}  }\right\vert =\left\vert \left\{
\left(  k,j\right)  \ \mid\ k\in\left\{  i+1,i+2,\ldots,j-1\right\}  \right\}
\right\vert =j-i-1,
\]
qed.}. Also, these sets $A$ and $B$ are disjoint\footnote{\textit{Proof.}
Assume the contrary. Then, $A\cap B\neq\varnothing$. Hence, there exists some
element $c$ of the set $A\cap B$. Consider such a $c$. We have%
\[
c\in A\cap B\subseteq A=\left\{  \left(  i,k\right)  \ \mid\ k\in\left\{
i+1,i+2,\ldots,j\right\}  \right\}  .
\]
In other words, we can write $c$ in the form $c=\left(  i,k\right)  $ for some
$k\in\left\{  i+1,i+2,\ldots,j\right\}  $. Let us denote this $k$ by $\ell$.
Then, $\ell$ is an element of $\left\{  i+1,i+2,\ldots,j\right\}  $ and
satisfies $c=\left(  i,\ell\right)  $.
\par
Also, $c\in A\cap B\subseteq B=\left\{  \left(  k,j\right)  \ \mid
\ k\in\left\{  i+1,i+2,\ldots,j-1\right\}  \right\}  $. Hence, we can write
$c$ in the form $c=\left(  k,j\right)  $ for some $k\in\left\{  i+1,i+2,\ldots
,j-1\right\}  $. Consider this $k$. Now, $c=\left(  k,j\right)  $, so that
$\left(  k,j\right)  =c=\left(  i,\ell\right)  $. Hence, $k=i$ and $j=\ell$.
Thus, $i=k\in\left\{  i+1,i+2,\ldots,j-1\right\}  $, which contradicts
$i\notin\left\{  i+1,i+2,\ldots,j-1\right\}  $. This contradiction shows that
our assumption was wrong, qed.}. Hence,
\[
\left\vert A\cup B\right\vert =\underbrace{\left\vert A\right\vert }%
_{=j-i}+\underbrace{\left\vert B\right\vert }_{=j-i-1}=\left(  j-i\right)
+\left(  j-i-1\right)  =2\underbrace{\left(  j-i\right)  }_{=\left\vert
j-i\right\vert }-1=2\left\vert j-i\right\vert -1.
\]

Now, let $\operatorname*{Inv}\left(  t_{i,j}\right)  $ denote the set of all
inversions of $t_{i,j}$. Then, $\ell\left(  t_{i,j}\right)  =\left\vert
\operatorname*{Inv}\left(  t_{i,j}\right)  \right\vert $%
\ \ \ \ \footnote{\textit{Proof.} Recall that $\ell\left(  t_{i,j}\right)  $
is defined as the number of inversions of $t_{i,j}$. Thus,%
\begin{align*}
\ell\left(  t_{i,j}\right)   &  =\left(  \text{the number of inversions of
}t_{i,j}\right) \\
&  =\left\vert \underbrace{\left(  \text{the set of all inversions of }%
t_{i,j}\right)  }_{=\operatorname*{Inv}\left(  t_{i,j}\right)  }\right\vert
=\left\vert \operatorname*{Inv}\left(  t_{i,j}\right)  \right\vert ,
\end{align*}
qed.}.

But $A\subseteq\operatorname*{Inv}\left(  t_{i,j}\right)  $%
\ \ \ \ \footnote{\textit{Proof.} Let $c\in A$. We shall show that
$c\in\operatorname*{Inv}\left(  t_{i,j}\right)  $.
\par
We have $c\in A=\left\{  \left(  i,k\right)  \ \mid\ k\in\left\{
i+1,i+2,\ldots,j\right\}  \right\}  $. Hence, $c$ can be written in the form
$c=\left(  i,k\right)  $ for some $k\in\left\{  i+1,i+2,\ldots,j\right\}  $.
Consider this $k$. From $k\in\left\{  i+1,i+2,\ldots,j\right\}  $, we obtain
$i<k\leq j$, so that $k\leq j\leq n$ and thus $1\leq i<k\leq n$.
\par
We shall now show that $j>t_{i,j}\left(  k\right)  $. In fact, we must be in
one of the following two cases:
\par
\textit{Case 1:} We have $k=j$.
\par
\textit{Case 2:} We have $k\neq j$.
\par
Let us first consider Case 1. In this case, we have $k=j$. Thus,
$t_{i,j}\left(  \underbrace{k}_{=j}\right)  =t_{i,j}\left(  j\right)  =i$.
Thus, $j>i=t_{i,j}\left(  k\right)  $. Hence, $j>t_{i,j}\left(  k\right)  $ is
proven in Case 1.
\par
Let us next consider Case 2. In this case, we have $k\neq j$. Combined with
$k\neq i$ (since $i<k$), this yields $k\notin\left\{  i,j\right\}  $. Hence,
$t_{i,j}\left(  k\right)  =k$ (by (\ref{sol.ps4.1ab.a.tij})). Also, combining
$k\leq j$ with $k\neq j$, we obtain $k<j$ and thus $j>k=t_{i,j}\left(
k\right)  $. Hence, $j>t_{i,j}\left(  k\right)  $ is proven in Case 2.
\par
Now, we have proven $j>t_{i,j}\left(  k\right)  $ in each of the two Cases 1
and 2. Thus, $j>t_{i,j}\left(  k\right)  $ always holds.
\par
Now, $t_{i,j}\left(  i\right)  =j>t_{i,j}\left(  k\right)  $. Hence, $\left(
i,k\right)  $ is a pair of integers satisfying $1\leq i<k\leq n$ and
$t_{i,j}\left(  i\right)  >t_{i,j}\left(  k\right)  $. In other words,
$\left(  i,k\right)  $ is an inversion of $t_{i,j}$ (by the definition of
\textquotedblleft inversion of $t_{i,j}$\textquotedblright). In other words,
$\left(  i,k\right)  \in\operatorname*{Inv}\left(  t_{i,j}\right)  $ (since
$\operatorname*{Inv}\left(  t_{i,j}\right)  $ is the set of all inversions of
$t_{i,j}$). Thus, $c=\left(  i,k\right)  \in\operatorname*{Inv}\left(
t_{i,j}\right)  $.
\par
Now, let us forget that we fixed $c$. We thus have shown that $c\in
\operatorname*{Inv}\left(  t_{i,j}\right)  $ for every $c\in A$. In other
words, $A\subseteq\operatorname*{Inv}\left(  t_{i,j}\right)  $, qed.} and
$B\subseteq\operatorname*{Inv}\left(  t_{i,j}\right)  $%
\ \ \ \ \footnote{\textit{Proof.} Let $c\in B$. We shall show that
$c\in\operatorname*{Inv}\left(  t_{i,j}\right)  $.
\par
We have $c\in B=\left\{  \left(  k,j\right)  \ \mid\ k\in\left\{
i+1,i+2,\ldots,j-1\right\}  \right\}  $. Hence, $c$ can be written in the form
$c=\left(  k,j\right)  $ for some $k\in\left\{  i+1,i+2,\ldots,j-1\right\}  $.
Consider this $k$. From $k\in\left\{  i+1,i+2,\ldots,j-1\right\}  $, we obtain
$i<k<j$, so that $1\leq i<k$ and thus $1\leq k<j\leq n$. Also, combining
$k\neq i$ (since $i<k$) and $k\neq j$ (since $k<j$), we obtain $k\notin%
\left\{  i,j\right\}  $. Hence, $t_{i,j}\left(  k\right)  =k$ (by
(\ref{sol.ps4.1ab.a.tij})). But $t_{i,j}\left(  j\right)  =i<k=t_{i,j}\left(
k\right)  $, so that $t_{i,j}\left(  k\right)  >t_{i,j}\left(  j\right)  $.
\par
Hence, $\left(  k,j\right)  $ is a pair of integers satisfying $1\leq k<j\leq
n$ and $t_{i,j}\left(  k\right)  >t_{i,j}\left(  j\right)  $. In other words,
$\left(  k,j\right)  $ is an inversion of $t_{i,j}$ (by the definition of
\textquotedblleft inversion of $t_{i,j}$\textquotedblright). In other words,
$\left(  k,j\right)  \in\operatorname*{Inv}\left(  t_{i,j}\right)  $ (since
$\operatorname*{Inv}\left(  t_{i,j}\right)  $ is the set of all inversions of
$t_{i,j}$). Thus, $c=\left(  k,j\right)  \in\operatorname*{Inv}\left(
t_{i,j}\right)  $.
\par
Now, let us forget that we fixed $c$. We thus have shown that $c\in
\operatorname*{Inv}\left(  t_{i,j}\right)  $ for every $c\in B$. In other
words, $B\subseteq\operatorname*{Inv}\left(  t_{i,j}\right)  $, qed.}. Hence,%
\[
\underbrace{A}_{\subseteq\operatorname*{Inv}\left(  t_{i,j}\right)  }%
\cup\underbrace{B}_{\subseteq\operatorname*{Inv}\left(  t_{i,j}\right)
}\subseteq\operatorname*{Inv}\left(  t_{i,j}\right)  \cup\operatorname*{Inv}%
\left(  t_{i,j}\right)  =\operatorname*{Inv}\left(  t_{i,j}\right)  .
\]

On the other hand, $\operatorname*{Inv}\left(  t_{i,j}\right)  \subseteq A\cup
B$\ \ \ \ \footnote{\textit{Proof.} Let $c\in\operatorname*{Inv}\left(
t_{i,j}\right)  $. We shall show that $c\in A\cup B$.
\par
We have $c\in\operatorname*{Inv}\left(  t_{i,j}\right)  $. In other words, $c$
is an inversion of $t_{i,j}$ (since $\operatorname*{Inv}\left(  t_{i,j}%
\right)  $ is the set of all inversions of $t_{i,j}$). In other words, $c$ is
a pair $\left(  u,v\right)  $ of integers satisfying $1\leq u<v\leq n$ and
$t_{i,j}\left(  u\right)  >t_{i,j}\left(  v\right)  $. Consider this $\left(
u,v\right)  $. Thus, $c=\left(  u,v\right)  $.
\par
We must be in one of the following two cases:
\par
\textit{Case 1:} We have $t_{i,j}\left(  u\right)  >u$.
\par
\textit{Case 2:} We don't have $t_{i,j}\left(  u\right)  >u$.
\par
Let us first consider Case 1. In this case, we have $t_{i,j}\left(  u\right)
>u$. If we had $u\notin\left\{  i,j\right\}  $, then we would have
$t_{i,j}\left(  u\right)  =u$ (by (\ref{sol.ps4.1ab.a.tij}), applied to
$k=u$), which would contradict $t_{i,j}\left(  u\right)  >u$. Hence, we cannot
have $u\notin\left\{  i,j\right\}  $. Thus, we have $u\in\left\{  i,j\right\}
$. If we had $u=j$, then we would have $t_{i,j}\left(  \underbrace{u}%
_{=j}\right)  =t_{i,j}\left(  j\right)  =i<j$, which would contradict
$t_{i,j}\left(  u\right)  >u=j$. Hence, we cannot have $u=j$. We thus have
$u\neq j$. Since $u\in\left\{  i,j\right\}  $ but $u\neq j$, we must have
$u=i$. Hence, $t_{i,j}\left(  \underbrace{u}_{=i}\right)  =t_{i,j}\left(
i\right)  =j$, so that $j=t_{i,j}\left(  u\right)  >t_{i,j}\left(  v\right)
$.
\par
Now, we assume (for the sake of contradiction) that $v>j$. Combining $v\neq i$
(since $v>j>i$) and $v\neq j$ (since $v>j$), we obtain $v\notin\left\{
i,j\right\}  $, so that $t_{i,j}\left(  v\right)  =v$ (by
(\ref{sol.ps4.1ab.a.tij}), applied to $k=v$). Hence, $j>t_{i,j}\left(
v\right)  =v>j$, which is absurd. Thus, we have obtained a contradiction.
Therefore, our assumption (that $v>j$) was wrong. We thus must have $v\leq j$.
Combined with $i=u<v$, this yields $i<v\leq j$, so that $v\in\left\{
i+1,i+2,\ldots,j\right\}  $. Now,%
\begin{align*}
c  &  =\left(  \underbrace{u}_{=i},v\right)  =\left(  i,v\right)  \in\left\{
\left(  i,k\right)  \ \mid\ k\in\left\{  i+1,i+2,\ldots,j\right\}  \right\} \\
&  \ \ \ \ \ \ \ \ \ \ \left(  \text{since }v\in\left\{  i+1,i+2,\ldots
,j\right\}  \right) \\
&  =A\subseteq A\cup B.
\end{align*}
Thus, $c\in A\cup B$ is proven in Case 1.
\par
Let us next consider Case 2. In this case, we don't have $t_{i,j}\left(
u\right)  >u$. Hence, we have $t_{i,j}\left(  u\right)  \leq u$, so that
$u\geq t_{i,j}\left(  u\right)  >t_{i,j}\left(  v\right)  $ and thus
$t_{i,j}\left(  v\right)  <u<v$. If we had $v\notin\left\{  i,j\right\}  $,
then we would have $t_{i,j}\left(  v\right)  =v$ (by (\ref{sol.ps4.1ab.a.tij}%
), applied to $k=v$), which would contradict $t_{i,j}\left(  v\right)  <v$.
Hence, we cannot have $v\notin\left\{  i,j\right\}  $. Thus, we have
$v\in\left\{  i,j\right\}  $. If we had $v=i$, then we would have
$t_{i,j}\left(  \underbrace{v}_{=i}\right)  =t_{i,j}\left(  i\right)  =j>i$,
which would contradict $t_{i,j}\left(  v\right)  <v=i$. Hence, we cannot have
$v=i$. We thus have $v\neq i$. Since $v\in\left\{  i,j\right\}  $ but $v\neq
i$, we must have $v=j$. Hence, $t_{i,j}\left(  \underbrace{v}_{=j}\right)
=t_{i,j}\left(  j\right)  =i$, so that $t_{i,j}\left(  u\right)
>t_{i,j}\left(  v\right)  =i$. Hence, $i<t_{i,j}\left(  u\right)  $.
\par
If we had $u=i$, then we would have $t_{i,j}\left(  \underbrace{u}%
_{=i}\right)  =t_{i,j}\left(  i\right)  =j>i=u$, which would contradict
$t_{i,j}\left(  u\right)  \leq u$. Thus, we cannot have $u=i$. Hence, $u\neq
i$.
\par
Now, we assume (for the sake of contradiction) that $u\leq i$. Combining
$u\neq i$ and $u\neq j$ (since $u\leq i<j$), we obtain $u\notin\left\{
i,j\right\}  $, so that $t_{i,j}\left(  u\right)  =u$ (by
(\ref{sol.ps4.1ab.a.tij}), applied to $k=u$). Hence, $i<t_{i,j}\left(
u\right)  =u\leq i$, which is absurd. Thus, we have obtained a contradiction.
Therefore, our assumption (that $u\leq i$) was wrong. We thus must have $u>i$.
In other words, $i<u$. Combined with $u<v\leq j$, this yields $i<u<j$, so that
$u\in\left\{  i+1,i+2,\ldots,j-1\right\}  $. Now,%
\begin{align*}
c  &  =\left(  u,\underbrace{v}_{=j}\right)  =\left(  u,j\right)  \in\left\{
\left(  k,j\right)  \ \mid\ k\in\left\{  i+1,i+2,\ldots,j-1\right\}  \right\}
\\
&  \ \ \ \ \ \ \ \ \ \ \left(  \text{since }v\in\left\{  i+1,i+2,\ldots
,j-1\right\}  \right) \\
&  =B\subseteq A\cup B.
\end{align*}
Thus, $c\in A\cup B$ is proven in Case 2.
\par
We therefore have shown that $c\in A\cup B$ in each of the two Cases 1 and 2.
Hence, $c\in A\cup B$ always holds.
\par
Now, let us forget that we fixed $c$. We thus have proven that $c\in A\cup B$
for every $c\in\operatorname*{Inv}\left(  t_{i,j}\right)  $. In other words,
$\operatorname*{Inv}\left(  t_{i,j}\right)  \subseteq A\cup B$, qed.}.
Combined with $A\cup B\subseteq\operatorname*{Inv}\left(  t_{i,j}\right)  $,
this yields $A\cup B=\operatorname*{Inv}\left(  t_{i,j}\right)  $. Thus,
$\left\vert A\cup B\right\vert =\left\vert \operatorname*{Inv}\left(
t_{i,j}\right)  \right\vert $. Compared with $\ell\left(  t_{i,j}\right)
=\left\vert \operatorname*{Inv}\left(  t_{i,j}\right)  \right\vert $, this
yields $\ell\left(  t_{i,j}\right)  =\left\vert A\cup B\right\vert
=2\left\vert j-i\right\vert -1$. This proves (\ref{sol.ps4.1ab.a.claim}).]

Exercise \ref{exe.ps4.1ab} \textbf{(a)} is thus solved.
\end{verlong}

\textbf{(b)} \textit{First solution to Exercise \ref{exe.ps4.1ab}
\textbf{(b)}:} From (\ref{sol.ps4.1ab.a.claim}), we have $\ell\left(
t_{i,j}\right)  =2\left\vert j-i\right\vert -1$.

But the integer $2\left\vert j-i\right\vert -1$ is odd. Thus, $\left(
-1\right)  ^{2\left\vert j-i\right\vert -1}=-1$. But the definition of
$\left(  -1\right)  ^{t_{i,j}}$ yields%
\begin{align*}
\left(  -1\right)  ^{t_{i,j}}  &  =\left(  -1\right)  ^{\ell\left(
t_{i,j}\right)  }=\left(  -1\right)  ^{2\left\vert j-i\right\vert
-1}\ \ \ \ \ \ \ \ \ \ \left(  \text{since }\ell\left(  t_{i,j}\right)
=2\left\vert j-i\right\vert -1\right) \\
&  =-1.
\end{align*}
This solves Exercise \ref{exe.ps4.1ab} \textbf{(b)}.

\textit{Second solution to Exercise \ref{exe.ps4.1ab} \textbf{(b)}:} Here is
an alternative solution of Exercise \ref{exe.ps4.1ab} \textbf{(b)} which makes
no use of part \textbf{(a)}.

The set $\left\{  1,2,\ldots,n\right\}  $ has at least two distinct elements
(namely, $i$ and $j$). Hence, $n\geq2$.

There exists a permutation $\sigma\in S_{n}$ such that $\left(  i,j\right)
=\left(  \sigma\left(  1\right)  ,\sigma\left(  2\right)  \right)
$\ \ \ \ \footnote{\textit{Proof.} Let $\left[  n\right]  =\left\{
1,2,\ldots,n\right\}  $. Notice that $n\geq2$ and thus $2\in\left\{
0,1,\ldots,n\right\}  $.
\par
The integers $i$ and $j$ are distinct. Hence, $\left(  i,j\right)  $ is a list
of some elements of $\left[  n\right]  $ such that $i$ and $j$ are distinct.
Therefore, Proposition \ref{prop.perms.lists} \textbf{(c)} (applied to $k=2$
and $\left(  p_{1},p_{2},\ldots,p_{k}\right)  =\left(  i,j\right)  $) yields
that there exists a permutation $\sigma\in S_{n}$ such that $\left(
i,j\right)  =\left(  \sigma\left(  1\right)  ,\sigma\left(  2\right)  \right)
$. Qed.}. Consider such a $\sigma$.

We have $\left(  i,j\right)  =\left(  \sigma\left(  1\right)  ,\sigma\left(
2\right)  \right)  $, thus $\left(  \sigma\left(  1\right)  ,\sigma\left(
2\right)  \right)  =\left(  i,j\right)  $. In other words, $\sigma\left(
1\right)  =i$ and $\sigma\left(  2\right)  =j$.

We have $n\geq2$. Hence, the permutation $s_{1}$ in $S_{n}$ is well-defined.
According to its definition, this permutation $s_{1}$ swaps $1$ with $2$ but
leaves all other numbers unchanged. In other words, we have $s_{1}\left(
1\right)  =2$, $s_{1}\left(  2\right)  =1$, and%
\begin{equation}
s_{1}\left(  k\right)  =k\ \ \ \ \ \ \ \ \ \ \text{for every }k\in\left\{
1,2,\ldots,n\right\}  \text{ satisfying }k\notin\left\{  1,2\right\}  .
\label{sol.ps4.1ab.b.s1}%
\end{equation}

On the other hand, the permutation $t_{i,j}$ swaps $i$ with $j$ while leaving
all other elements of $\left\{  1,2,\ldots,n\right\}  $ unchanged. In other
words, we have $t_{i,j}\left(  i\right)  =j$, $t_{i,j}\left(  j\right)  =i$
and%
\begin{equation}
t_{i,j}\left(  k\right)  =k\ \ \ \ \ \ \ \ \ \ \text{for every }k\in\left\{
1,2,\ldots,n\right\}  \text{ satisfying }k\notin\left\{  i,j\right\}  .
\label{sol.ps4.1ab.b.tij}%
\end{equation}

\begin{vershort}
Now, every $k\in\left\{  1,2,\ldots,n\right\}  $ satisfies $t_{i,j}\left(
\sigma\left(  k\right)  \right)  =\sigma\left(  s_{1}\left(  k\right)
\right)  $\ \ \ \ \footnote{\textit{Proof.} Let $k\in\left\{  1,2,\ldots
,n\right\}  $. We need to show that $t_{i,j}\left(  \sigma\left(  k\right)
\right)  =\sigma\left(  s_{1}\left(  k\right)  \right)  $.
\par
We are in one of the following three cases:
\par
\textit{Case 1:} We have $k=1$.
\par
\textit{Case 2:} We have $k=2$.
\par
\textit{Case 3:} We have $k\notin\left\{  1,2\right\}  $.
\par
Let us first consider Case 1. In this case, we have $k=1$. Hence,
$t_{i,j}\left(  \sigma\left(  \underbrace{k}_{=1}\right)  \right)
=t_{i,j}\left(  \underbrace{\sigma\left(  1\right)  }_{=i}\right)
=t_{i,j}\left(  i\right)  =j$. Compared with $\sigma\left(  s_{1}\left(
\underbrace{k}_{=1}\right)  \right)  =\sigma\left(  \underbrace{s_{1}\left(
1\right)  }_{=2}\right)  =\sigma\left(  2\right)  =j$, this yields
$t_{i,j}\left(  \sigma\left(  k\right)  \right)  =\sigma\left(  s_{1}\left(
k\right)  \right)  $. Hence, $t_{i,j}\left(  \sigma\left(  k\right)  \right)
=\sigma\left(  s_{1}\left(  k\right)  \right)  $ is proven in Case 1.
\par
The proof of $t_{i,j}\left(  \sigma\left(  k\right)  \right)  =\sigma\left(
s_{1}\left(  k\right)  \right)  $ in Case 2 is similar and left to the reader.
\par
Let us first consider Case 3. In this case, we have $k\notin\left\{
1,2\right\}  $. Hence, $s_{1}\left(  k\right)  =k$ (by (\ref{sol.ps4.1ab.b.s1}%
)). On the other hand, the map $\sigma$ is a permutation (since $\sigma\in
S_{n}$), thus injective. Now, from $k\neq1$ (since $k\notin\left\{
1,2\right\}  $), we obtain $\sigma\left(  k\right)  \neq\sigma\left(
1\right)  $ (since the map $\sigma$ is injective), so that $\sigma\left(
k\right)  \neq\sigma\left(  1\right)  =i$. Similarly, from $k\neq2$, we can
obtain $\sigma\left(  k\right)  \neq j$. Combining $\sigma\left(  k\right)
\neq i$ with $\sigma\left(  k\right)  \neq j$, we obtain $\sigma\left(
k\right)  \notin\left\{  i,j\right\}  $, and therefore $t_{i,j}\left(
\sigma\left(  k\right)  \right)  =\sigma\left(  k\right)  $ (by
(\ref{sol.ps4.1ab.b.tij}), applied to $\sigma\left(  k\right)  $ instead of
$k$). Compared with $\sigma\left(  \underbrace{s_{1}\left(  k\right)  }%
_{=k}\right)  =\sigma\left(  k\right)  $, this yields $t_{i,j}\left(
\sigma\left(  k\right)  \right)  =\sigma\left(  s_{1}\left(  k\right)
\right)  $. Hence, $t_{i,j}\left(  \sigma\left(  k\right)  \right)
=\sigma\left(  s_{1}\left(  k\right)  \right)  $ is proven in Case 3.
\par
Thus, $t_{i,j}\left(  \sigma\left(  k\right)  \right)  =\sigma\left(
s_{1}\left(  k\right)  \right)  $ is proven in each of the three Cases 1, 2
and 3. Hence, $t_{i,j}\left(  \sigma\left(  k\right)  \right)  =\sigma\left(
s_{1}\left(  k\right)  \right)  $ always holds, qed.}. Hence, every
$k\in\left\{  1,2,\ldots,n\right\}  $ satisfies $\left(  t_{i,j}\circ
\sigma\right)  \left(  k\right)  =t_{i,j}\left(  \sigma\left(  k\right)
\right)  =\sigma\left(  s_{1}\left(  k\right)  \right)  =\left(  \sigma\circ
s_{1}\right)  \left(  k\right)  $. In other words, $t_{i,j}\circ\sigma
=\sigma\circ s_{1}$.
\end{vershort}

\begin{verlong}
Now, every $k\in\left\{  1,2,\ldots,n\right\}  $ satisfies $t_{i,j}\left(
\sigma\left(  k\right)  \right)  =\sigma\left(  s_{1}\left(  k\right)
\right)  $\ \ \ \ \footnote{\textit{Proof.} Let $k\in\left\{  1,2,\ldots
,n\right\}  $. We need to show that $t_{i,j}\left(  \sigma\left(  k\right)
\right)  =\sigma\left(  s_{1}\left(  k\right)  \right)  $.
\par
We are in one of the following three cases:
\par
\textit{Case 1:} We have $k=1$.
\par
\textit{Case 2:} We have $k=2$.
\par
\textit{Case 3:} We have $k\notin\left\{  1,2\right\}  $.
\par
Let us first consider Case 1. In this case, we have $k=1$. Hence,
$t_{i,j}\left(  \sigma\left(  \underbrace{k}_{=1}\right)  \right)
=t_{i,j}\left(  \underbrace{\sigma\left(  1\right)  }_{=i}\right)
=t_{i,j}\left(  i\right)  =j$. Compared with $\sigma\left(  s_{1}\left(
\underbrace{k}_{=1}\right)  \right)  =\sigma\left(  \underbrace{s_{1}\left(
1\right)  }_{=2}\right)  =\sigma\left(  2\right)  =j$, this yields
$t_{i,j}\left(  \sigma\left(  k\right)  \right)  =\sigma\left(  s_{1}\left(
k\right)  \right)  $. Hence, $t_{i,j}\left(  \sigma\left(  k\right)  \right)
=\sigma\left(  s_{1}\left(  k\right)  \right)  $ is proven in Case 1.
\par
Let us next consider Case 2. In this case, we have $k=2$. Hence,
$t_{i,j}\left(  \sigma\left(  \underbrace{k}_{=2}\right)  \right)
=t_{i,j}\left(  \underbrace{\sigma\left(  2\right)  }_{=j}\right)
=t_{i,j}\left(  j\right)  =i$. Compared with $\sigma\left(  s_{1}\left(
\underbrace{k}_{=2}\right)  \right)  =\sigma\left(  \underbrace{s_{1}\left(
2\right)  }_{=1}\right)  =\sigma\left(  1\right)  =i$, this yields
$t_{i,j}\left(  \sigma\left(  k\right)  \right)  =\sigma\left(  s_{1}\left(
k\right)  \right)  $. Hence, $t_{i,j}\left(  \sigma\left(  k\right)  \right)
=\sigma\left(  s_{1}\left(  k\right)  \right)  $ is proven in Case 2.
\par
Let us first consider Case 3. In this case, we have $k\notin\left\{
1,2\right\}  $. Hence, $s_{1}\left(  k\right)  =k$ (by (\ref{sol.ps4.1ab.b.s1}%
)). On the other hand, the map $\sigma$ is a permutation (since $\sigma\in
S_{n}$), thus injective. Now, from $k\neq1$ (since $k\notin\left\{
1,2\right\}  $), we obtain $\sigma\left(  k\right)  \neq\sigma\left(
1\right)  $ (since the map $\sigma$ is injective), so that $\sigma\left(
k\right)  \neq\sigma\left(  1\right)  =i$. Also, from $k\neq2$ (since
$k\notin\left\{  1,2\right\}  $), we obtain $\sigma\left(  k\right)
\neq\sigma\left(  2\right)  $ (since the map $\sigma$ is injective), so that
$\sigma\left(  k\right)  \neq\sigma\left(  2\right)  =j$. Combining
$\sigma\left(  k\right)  \neq i$ with $\sigma\left(  k\right)  \neq j$, we
obtain $\sigma\left(  k\right)  \notin\left\{  i,j\right\}  $, and therefore
$t_{i,j}\left(  \sigma\left(  k\right)  \right)  =\sigma\left(  k\right)  $
(by (\ref{sol.ps4.1ab.b.tij}), applied to $\sigma\left(  k\right)  $ instead
of $k$). Compared with $\sigma\left(  \underbrace{s_{1}\left(  k\right)
}_{=k}\right)  =\sigma\left(  k\right)  $, this yields $t_{i,j}\left(
\sigma\left(  k\right)  \right)  =\sigma\left(  s_{1}\left(  k\right)
\right)  $. Hence, $t_{i,j}\left(  \sigma\left(  k\right)  \right)
=\sigma\left(  s_{1}\left(  k\right)  \right)  $ is proven in Case 3.
\par
Thus, $t_{i,j}\left(  \sigma\left(  k\right)  \right)  =\sigma\left(
s_{1}\left(  k\right)  \right)  $ is proven in each of the three Cases 1, 2
and 3. Hence, $t_{i,j}\left(  \sigma\left(  k\right)  \right)  =\sigma\left(
s_{1}\left(  k\right)  \right)  $ always holds, qed.}. Hence, every
$k\in\left\{  1,2,\ldots,n\right\}  $ satisfies $\left(  t_{i,j}\circ
\sigma\right)  \left(  k\right)  =t_{i,j}\left(  \sigma\left(  k\right)
\right)  =\sigma\left(  s_{1}\left(  k\right)  \right)  =\left(  \sigma\circ
s_{1}\right)  \left(  k\right)  $. In other words, $t_{i,j}\circ\sigma
=\sigma\circ s_{1}$.
\end{verlong}

Proposition \ref{prop.perm.signs.basics} \textbf{(b)} shows that $\left(
-1\right)  ^{s_{k}}=-1$ for every $k\in\left\{  1,2,\ldots,n-1\right\}  $.
Applied to $k=1$, this yields $\left(  -1\right)  ^{s_{1}}=-1$.

On the other hand, (\ref{eq.sign.prod}) (applied to $\tau=s_{1}$) yields
$\left(  -1\right)  ^{\sigma\circ s_{1}}=\left(  -1\right)  ^{\sigma}%
\cdot\underbrace{\left(  -1\right)  ^{s_{1}}}_{=-1}=\left(  -1\right)
^{\sigma}\cdot\left(  -1\right)  $.

But (\ref{eq.sign.prod}) (applied to $t_{i,j}$ and $\sigma$ instead of
$\sigma$ and $\tau$) yields $\left(  -1\right)  ^{t_{i,j}\circ\sigma}=\left(
-1\right)  ^{t_{i,j}}\cdot\left(  -1\right)  ^{\sigma}$, so that%
\begin{align*}
\left(  -1\right)  ^{t_{i,j}}\cdot\left(  -1\right)  ^{\sigma}  &  =\left(
-1\right)  ^{t_{i,j}\circ\sigma}=\left(  -1\right)  ^{\sigma\circ s_{1}%
}\ \ \ \ \ \ \ \ \ \ \left(  \text{since }t_{i,j}\circ\sigma=\sigma\circ
s_{1}\right) \\
&  =\left(  -1\right)  ^{\sigma}\cdot\left(  -1\right)  .
\end{align*}
We can cancel $\left(  -1\right)  ^{\sigma}$ from this equality (since
$\left(  -1\right)  ^{\sigma}\in\left\{  1,-1\right\}  $ is a nonzero
integer), and thus obtain $\left(  -1\right)  ^{t_{i,j}}=-1$. This solves
Exercise \ref{exe.ps4.1ab} \textbf{(b)} again.
\end{proof}

\subsection{Solution to Exercise \ref{exe.ps4.1c}}

\begin{vershort}
\begin{proof}
[Solution to Exercise \ref{exe.ps4.1c}.]The inversions of $w_{0}$ are the
pairs $\left(  i,j\right)  $ of integers satisfying $1\leq i<j\leq n$ and
$w_{0}\left(  i\right)  >w_{0}\left(  j\right)  $ (because this is how we
defined inversions). Since \textbf{every} pair of integers $\left(
i,j\right)  $ satisfying $1\leq i<j\leq n$ automatically satisfies
$w_{0}\left(  i\right)  >w_{0}\left(  j\right)  $%
\ \ \ \ \footnote{\textit{Proof.} Let $\left(  i,j\right)  $ be a pair of
integers satisfying $1\leq i<j\leq n$. The definition of $w_{0}$ yields
$w_{0}\left(  i\right)  =n+1-i$ and $w_{0}\left(  j\right)  =n+1-j$. Hence,
$w_{0}\left(  i\right)  =n+1-\underbrace{i}_{<j}>n+1-j=w_{0}\left(  j\right)
$, qed.}, we can simplify this statement as follows: The inversions of $w_{0}$
are the pairs $\left(  i,j\right)  $ of integers satisfying $1\leq i<j\leq n$.
But the number of such pairs is $n\left(  n-1\right)  /2$%
\ \ \ \ \footnote{There are two ways to prove this:
\par
\begin{itemize}
\item Either we can argue that these pairs are in a one-to-one correspondence
with the $2$-element subsets of $\left\{  1,2,\ldots,n\right\}  $. (Namely,
any pair $\left(  i,j\right)  $ corresponds to the subset $\left\{
i,j\right\}  $, and conversely, any subset $S$ corresponds to the pair
$\left(  \min S,\max S\right)  $.) Therefore, the number of such pairs equals
the number of all $2$-element subsets of $\left\{  1,2,\ldots,n\right\}  $;
but the latter number is known to be $\dbinom{n}{2}=n\left(  n-1\right)  /2$.
\par
\item Alternatively, we can compute this number as follows: For every pair
$\left(  i,j\right)  $ of integers satisfying $1\leq i<j\leq n$, we have
$j\in\left\{  2,3,\ldots,n\right\}  $ (since $1<j\leq n$). Hence,
\begin{align*}
&  \left(  \text{the number of pairs }\left(  i,j\right)  \text{ of integers
satisfying }1\leq i<j\leq n\right) \\
&  =\sum_{k=2}^{n}\underbrace{\left(  \text{the number of pairs } \left(  i,
j\right)  \text{ satisfying }1\leq i<j\leq n\text{ and } j=k\right)
}_{=\left(  \text{the number of integers }i\text{ satisfying }1\leq
i<k\right)  = k-1 }\\
&  =\sum_{k=2}^{n}\left(  k-1\right)  =\sum_{j=1}^{n-1}%
j\ \ \ \ \ \ \ \ \ \ \left(  \text{here, we substituted }j\text{ for
}k-1\text{ in the sum}\right) \\
&  = \sum_{i=1}^{n-1} i = \dfrac{\left(  n-1\right)  \left(  \left(
n-1\right)  +1\right)  }{2}\\
&  \ \ \ \ \ \ \ \ \ \ \left(  \text{by (\ref{eq.sum.littlegauss2}), applied
to $n-1$ instead of $n$}\right) \\
&  =n\left(  n-1\right)  /2.
\end{align*}
\end{itemize}
}. Thus, the number of inversions of $w_{0}$ is $n\left(  n-1\right)  /2$. In
other words, $\ell\left(  w_{0}\right)  =n\left(  n-1\right)  /2$ (since
$\ell\left(  w_{0}\right)  $ is defined as the number of inversions of $w_{0}%
$). Therefore, the definition of $\left(  -1\right)  ^{w_{0}}$ yields $\left(
-1\right)  ^{w_{0}}=\left(  -1\right)  ^{\ell\left(  w_{0}\right)  }=\left(
-1\right)  ^{n\left(  n-1\right)  /2}$ (since $\ell\left(  w_{0}\right)
=n\left(  n-1\right)  /2$).

At this point, we could declare Exercise \ref{exe.ps4.1c} to be solved, since
we have found formulas for both $\ell\left(  w_{0}\right)  $ and $\left(
-1\right)  ^{w_{0}}$. Nevertheless, let us give a different expression for
$\left(  -1\right)  ^{w_{0}}$, which can be evaluated faster. Namely, we claim
that%
\begin{equation}
\left(  -1\right)  ^{w_{0}}=
\begin{cases}
1, & \text{if }n\equiv0\operatorname{mod}4\text{ or }n\equiv
1\operatorname{mod}4;\\
-1, & \text{if }n\equiv2\operatorname{mod}4\text{ or }n\equiv
3\operatorname{mod}4
\end{cases}
. \label{sol.ps4.1c.short.extra}%
\end{equation}

[\textit{Proof of (\ref{sol.ps4.1c.short.extra}):} We must be in one of the
following four cases:

\textit{Case 1:} We have $n\equiv0\operatorname{mod}4$.

\textit{Case 2:} We have $n\equiv1\operatorname{mod}4$.

\textit{Case 3:} We have $n\equiv2\operatorname{mod}4$.

\textit{Case 4:} We have $n\equiv3\operatorname{mod}4$.

The proofs of (\ref{sol.ps4.1c.short.extra}) in these four cases are more or
less analogous. Let us only show the proof in Case 4. In this case, we have
$n\equiv3\operatorname{mod}4$. Thus, $n=4m+3$ for some $m\in\mathbb{Z}$.
Consider this $m$. We have%
\begin{align*}
\underbrace{n}_{=4m+3}\underbrace{\left(  n-1\right)  }%
_{\substack{=4m+2\\\text{(since }n=4m+3\text{)}}}/2  &  =\left(  4m+3\right)
\underbrace{\left(  4m+2\right)  /2}_{=2m+1}=\left(  4m+3\right)  \left(
2m+1\right) \\
&  =8m^{2}+10m+3=2\left(  4m^{2}+5m+1\right)  +1.
\end{align*}
Thus, the integer $n\left(  n-1\right)  /2$ is odd, so that $\left(
-1\right)  ^{n\left(  n-1\right)  /2}=-1$. Hence,%
\[
\left(  -1\right)  ^{w_{0}}=\left(  -1\right)  ^{n\left(  n-1\right)  /2}=-1.
\]
Compared with $%
\begin{cases}
1, & \text{if }n\equiv0\operatorname{mod}4\text{ or }n\equiv
1\operatorname{mod}4;\\
-1, & \text{if }n\equiv2\operatorname{mod}4\text{ or }n\equiv
3\operatorname{mod}4
\end{cases}
=-1$ (since $n\equiv3\operatorname{mod}4$), this yields $\left(  -1\right)
^{w_{0}}=%
\begin{cases}
1, & \text{if }n\equiv0\operatorname{mod}4\text{ or }n\equiv
1\operatorname{mod}4;\\
-1, & \text{if }n\equiv2\operatorname{mod}4\text{ or }n\equiv
3\operatorname{mod}4
\end{cases}
$. Thus, (\ref{sol.ps4.1c.short.extra}) is proven in Case 4. The other three
cases are analogous, and so we conclude that (\ref{sol.ps4.1c.short.extra}) holds.]
\end{proof}
\end{vershort}

\begin{verlong}
\begin{proof}
[Solution to Exercise \ref{exe.ps4.1c}.]Let $G=\left\{  \left(  i,j\right)
\in\mathbb{Z}^{2}\ \mid\ 1\leq i<j\leq n\right\}  $. Then, $\left\vert
G\right\vert =n\left(  n-1\right)  /2$\ \ \ \ \footnote{\textit{Proof.} This
is fairly obvious, but let us nevertheless give a proof for the sake of
completeness.
\par
From (\ref{eq.sum.littlegauss2}), we obtain $\sum_{i=1}^{n}i=\dfrac{n\left(
n+1\right)  }{2}=n\left(  n+1\right)  /2$.
\par
We have%
\[
\sum_{\substack{\left(  i,j\right)  \in\mathbb{Z}^{2};\\1\leq i<j\leq
n}}1=\left\vert \underbrace{\left\{  \left(  i,j\right)  \in\mathbb{Z}%
^{2}\ \mid\ 1\leq i<j\leq n\right\}  }_{=G}\right\vert \cdot1=\left\vert
G\right\vert \cdot1=\left\vert G\right\vert ,
\]
so that%
\begin{align*}
\left\vert G\right\vert  &  =\underbrace{\sum_{\substack{\left(  i,j\right)
\in\mathbb{Z}^{2};\\1\leq i<j\leq n}}}_{\substack{=\sum_{\substack{\left(
i,j\right)  \in\mathbb{Z}^{2};\\1\leq i\leq n;\ i<j\leq n}}\\\text{(since the
condition }1\leq i<j\leq n\\\text{is equivalent to }\left(  1\leq i\leq
n\text{ and }i<j\leq n\right)  \text{)}}}1=\underbrace{\sum_{\substack{\left(
i,j\right)  \in\mathbb{Z}^{2};\\1\leq i\leq n;\ i<j\leq n}}}_{=\sum
_{\substack{i\in\mathbb{Z};\\1\leq i\leq n}}\sum_{\substack{j\in
\mathbb{Z};\\i<j\leq n}}}1=\sum_{\substack{i\in\mathbb{Z};\\1\leq i\leq
n}}\underbrace{\sum_{\substack{j\in\mathbb{Z};\\i<j\leq n}}1}%
_{\substack{=\left\vert \left\{  j\in\mathbb{Z}\ \mid\ i<j\leq n\right\}
\right\vert \cdot1\\=\left\vert \left\{  j\in\mathbb{Z}\ \mid\ i<j\leq
n\right\}  \right\vert }}\\
&  =\sum_{\substack{i\in\mathbb{Z};\\1\leq i\leq n}}\left\vert
\underbrace{\left\{  j\in\mathbb{Z}\ \mid\ i<j\leq n\right\}  }_{=\left\{
i+1,i+2,\ldots,n\right\}  }\right\vert =\underbrace{\sum_{\substack{i\in
\mathbb{Z};\\1\leq i\leq n}}}_{=\sum_{i=1}^{n}}\underbrace{\left\vert \left\{
i+1,i+2,\ldots,n\right\}  \right\vert }_{\substack{=n-i\\\text{(since }i\leq
n\text{)}}}\\
&  =\sum_{i=1}^{n}\left(  n-i\right)  =\underbrace{\sum_{i=1}^{n}n}%
_{=nn}-\underbrace{\sum_{i=1}^{n}i}_{=n\left(  n+1\right)  /2}=nn-n\left(
n+1\right)  /2=n\left(  n-1\right)  /2,
\end{align*}
qed.}.

Now, let $\operatorname*{Inv}\left(  w_{0}\right)  $ denote the set of all
inversions of $w_{0}$. Then, $\ell\left(  w_{0}\right)  =\left\vert
\operatorname*{Inv}\left(  w_{0}\right)  \right\vert $%
\ \ \ \ \footnote{\textit{Proof.} Recall that $\ell\left(  w_{0}\right)  $ is
defined as the number of inversions of $w_{0}$. Thus,%
\begin{align*}
\ell\left(  w_{0}\right)   &  =\left(  \text{the number of inversions of
}w_{0}\right) \\
&  =\left\vert \underbrace{\left(  \text{the set of all inversions of }%
w_{0}\right)  }_{=\operatorname*{Inv}\left(  w_{0}\right)  }\right\vert
=\left\vert \operatorname*{Inv}\left(  w_{0}\right)  \right\vert ,
\end{align*}
qed.}.

On the other hand, $\operatorname*{Inv}\left(  w_{0}\right)  \subseteq
G$\ \ \ \ \footnote{\textit{Proof.} Let $c\in\operatorname*{Inv}\left(
w_{0}\right)  $. We shall show that $c\in G$.
\par
We have $c\in\operatorname*{Inv}\left(  w_{0}\right)  $. In other words, $c$
is an inversion of $w_{0}$ (since $\operatorname*{Inv}\left(  w_{0}\right)  $
is the set of all inversions of $w_{0}$). In other words, $c$ is a pair
$\left(  u,v\right)  $ of integers satisfying $1\leq u<v\leq n$ and
$w_{0}\left(  u\right)  >w_{0}\left(  v\right)  $. Consider this $\left(
u,v\right)  $. We have $c=\left(  u,v\right)  $ and $1\leq u<v\leq n$. Thus,
$c$ has the form $c=\left(  i,j\right)  $ for a pair $\left(  i,j\right)
\in\mathbb{Z}^{2}$ satisfying $1\leq i<j\leq n$ (namely, $\left(  i,j\right)
=\left(  u,v\right)  $). In other words, $c\in\left\{  \left(  i,j\right)
\in\mathbb{Z}^{2}\ \mid\ 1\leq i<j\leq n\right\}  =G$.
\par
Let us now forget that we fixed $c$. We thus have shown that $c\in G$ for
every $c\in\operatorname*{Inv}\left(  w_{0}\right)  $. In other words,
$\operatorname*{Inv}\left(  w_{0}\right)  \subseteq G$, qed.} and
$G\subseteq\operatorname*{Inv}\left(  w_{0}\right)  $%
\ \ \ \ \footnote{\textit{Proof.} Let $c\in G$. We shall prove that
$c\in\operatorname*{Inv}\left(  w_{0}\right)  $.
\par
We have $c\in G=\left\{  \left(  i,j\right)  \in\mathbb{Z}^{2}\ \mid\ 1\leq
i<j\leq n\right\}  $. In other words, $c$ can be written in the form
$c=\left(  i,j\right)  $ for some $\left(  i,j\right)  \in\mathbb{Z}^{2}$
satisfying $1\leq i<j\leq n$. Consider this $\left(  i,j\right)  $.
\par
The definition of $w_{0}$ yields $w_{0}\left(  j\right)  =n+1-j$ and
$w_{0}\left(  i\right)  =n+1-\underbrace{i}_{<j}>n+1-j=w_{0}\left(  j\right)
$. Hence, $\left(  i,j\right)  $ is a pair of integers satisfying $1\leq
i<j\leq n$ and $w_{0}\left(  i\right)  >w_{0}\left(  j\right)  $. In other
words, $\left(  i,j\right)  $ is an inversion of $w_{0}$ (by the definition of
an \textquotedblleft inversion of $w_{0}$\textquotedblright). In other words,
$\left(  i,j\right)  \in\operatorname*{Inv}\left(  w_{0}\right)  $ (since
$\operatorname*{Inv}\left(  w_{0}\right)  $ is the set of all inversions of
$w_{0}$). Thus, $c=\left(  i,j\right)  \in\operatorname*{Inv}\left(
w_{0}\right)  $.
\par
Now, let us forget that we fixed $c$. We thus have shown that $c\in
\operatorname*{Inv}\left(  w_{0}\right)  $ for every $c\in G$. In other words,
$G\subseteq\operatorname*{Inv}\left(  w_{0}\right)  $, qed.}. Combining these
two relations, we obtain $\operatorname*{Inv}\left(  w_{0}\right)  =G$. Thus,
$\left\vert \underbrace{\operatorname*{Inv}\left(  w_{0}\right)  }%
_{=G}\right\vert =\left\vert G\right\vert =n\left(  n-1\right)  /2$, so that
$\ell\left(  w_{0}\right)  =\left\vert \operatorname*{Inv}\left(
w_{0}\right)  \right\vert =n\left(  n-1\right)  /2$.

Now, the definition of $\left(  -1\right)  ^{w_{0}}$ yields $\left(
-1\right)  ^{w_{0}}=\left(  -1\right)  ^{\ell\left(  w_{0}\right)  }=\left(
-1\right)  ^{n\left(  n-1\right)  /2}$ (since $\ell\left(  w_{0}\right)
=n\left(  n-1\right)  /2$).

At this point, we could declare Exercise \ref{exe.ps4.1c} to be solved, since
we have found formulas for both $\ell\left(  w_{0}\right)  $ and $\left(
-1\right)  ^{w_{0}}$. Nevertheless, let us give a different expression for
$\left(  -1\right)  ^{w_{0}}$, which can be evaluated faster. Namely, we claim
that%
\begin{equation}
\left(  -1\right)  ^{w_{0}}=
\begin{cases}
1, & \text{if }n\equiv0\operatorname{mod}4\text{ or }n\equiv
1\operatorname{mod}4;\\
-1, & \text{if }n\equiv2\operatorname{mod}4\text{ or }n\equiv
3\operatorname{mod}4
\end{cases}
. \label{sol.ps4.1c.extra}%
\end{equation}

[\textit{Proof of (\ref{sol.ps4.1c.extra}):} We must be in one of the
following four cases:

\textit{Case 1:} We have $n\equiv0\operatorname{mod}4$.

\textit{Case 2:} We have $n\equiv1\operatorname{mod}4$.

\textit{Case 3:} We have $n\equiv2\operatorname{mod}4$.

\textit{Case 4:} We have $n\equiv3\operatorname{mod}4$.

Let us first consider Case 1. In this case, we have $n\equiv
0\operatorname{mod}4$. Thus, $n=4m$ for some $m\in\mathbb{Z}$. Consider this
$m$. We have $\underbrace{n}_{=4m}\left(  n-1\right)  /2=4m\left(  n-1\right)
/2=2m\left(  n-1\right)  $. Thus, the integer $n\left(  n-1\right)  /2$ is
even, so that $\left(  -1\right)  ^{n\left(  n-1\right)  /2}=1$. Hence,%
\[
\left(  -1\right)  ^{w_{0}}=\left(  -1\right)  ^{n\left(  n-1\right)  /2}=1.
\]
Compared with $%
\begin{cases}
1, & \text{if }n\equiv0\operatorname{mod}4\text{ or }n\equiv
1\operatorname{mod}4;\\
-1, & \text{if }n\equiv2\operatorname{mod}4\text{ or }n\equiv
3\operatorname{mod}4
\end{cases}
=1$ (since $n\equiv0\operatorname{mod}4$ or $n\equiv1\operatorname{mod}4$
(since $n\equiv0\operatorname{mod}4$)), this yields
\[
\left(  -1\right)  ^{w_{0}}=
\begin{cases}
1, & \text{if }n\equiv0\operatorname{mod}4\text{ or }n\equiv
1\operatorname{mod}4;\\
-1, & \text{if }n\equiv2\operatorname{mod}4\text{ or }n\equiv
3\operatorname{mod}4
\end{cases}
.
\]
Thus, (\ref{sol.ps4.1c.extra}) is proven in Case 1.

Let us next consider Case 2. In this case, we have $n\equiv1\operatorname{mod}%
4$. Thus, $n=4m+1$ for some $m\in\mathbb{Z}$. Consider this $m$. We have
$n\underbrace{\left(  n-1\right)  }_{\substack{=4m\\\text{(since
}n=4m+1\text{)}}}/2=n\cdot4m/2=2nm$. Thus, the integer $n\left(  n-1\right)
/2$ is even, so that $\left(  -1\right)  ^{n\left(  n-1\right)  /2}=1$. Hence,%
\[
\left(  -1\right)  ^{w_{0}}=\left(  -1\right)  ^{n\left(  n-1\right)  /2}=1.
\]
Compared with $%
\begin{cases}
1, & \text{if }n\equiv0\operatorname{mod}4\text{ or }n\equiv
1\operatorname{mod}4;\\
-1, & \text{if }n\equiv2\operatorname{mod}4\text{ or }n\equiv
3\operatorname{mod}4
\end{cases}
=1$ (since $n\equiv0\operatorname{mod}4$ or $n\equiv1\operatorname{mod}4$
(since $n\equiv1\operatorname{mod}4$)), this yields
\[
\left(  -1\right)  ^{w_{0}}=
\begin{cases}
1, & \text{if }n\equiv0\operatorname{mod}4\text{ or }n\equiv
1\operatorname{mod}4;\\
-1, & \text{if }n\equiv2\operatorname{mod}4\text{ or }n\equiv
3\operatorname{mod}4
\end{cases}
.
\]
Thus, (\ref{sol.ps4.1c.extra}) is proven in Case 2.

Let us next consider Case 3. In this case, we have $n\equiv2\operatorname{mod}%
4$. Thus, $n=4m+2$ for some $m\in\mathbb{Z}$. Consider this $m$. We have%
\begin{align*}
\underbrace{n}_{\substack{=4m+2\\=2\left(  2m+1\right)  }}\underbrace{\left(
n-1\right)  }_{\substack{=4m+1\\\text{(since }n=4m+2\text{)}}}/2  &  =2\left(
2m+1\right)  \left(  4m+1\right)  /2=\left(  2m+1\right)  \left(  4m+1\right)
\\
&  =8m^{2}+6m+1=2\left(  4m^{2}+3m\right)  +1.
\end{align*}
Thus, the integer $n\left(  n-1\right)  /2$ is odd, so that $\left(
-1\right)  ^{n\left(  n-1\right)  /2}=-1$. Hence,%
\[
\left(  -1\right)  ^{w_{0}}=\left(  -1\right)  ^{n\left(  n-1\right)  /2}=-1.
\]
Compared with $%
\begin{cases}
1, & \text{if }n\equiv0\operatorname{mod}4\text{ or }n\equiv
1\operatorname{mod}4;\\
-1, & \text{if }n\equiv2\operatorname{mod}4\text{ or }n\equiv
3\operatorname{mod}4
\end{cases}
=-1$ (since $n\equiv2\operatorname{mod}4$ or $n\equiv3\operatorname{mod}4$
(since $n\equiv2\operatorname{mod}4$)), this yields
\[
\left(  -1\right)  ^{w_{0}}=
\begin{cases}
1, & \text{if }n\equiv0\operatorname{mod}4\text{ or }n\equiv
1\operatorname{mod}4;\\
-1, & \text{if }n\equiv2\operatorname{mod}4\text{ or }n\equiv
3\operatorname{mod}4
\end{cases}
.
\]
Thus, (\ref{sol.ps4.1c.extra}) is proven in Case 3.

Let us finally consider Case 4. In this case, we have $n\equiv
3\operatorname{mod}4$. Thus, $n=4m+3$ for some $m\in\mathbb{Z}$. Consider this
$m$. We have%
\begin{align*}
\underbrace{n}_{=4m+3}\underbrace{\left(  n-1\right)  }%
_{\substack{=4m+2\\\text{(since }n=4m+3\text{)}}}/2  &  =\left(  4m+3\right)
\underbrace{\left(  4m+2\right)  /2}_{=2m+1}=\left(  4m+3\right)  \left(
2m+1\right) \\
&  =8m^{2}+10m+3=2\left(  4m^{2}+5m+1\right)  +1.
\end{align*}
Thus, the integer $n\left(  n-1\right)  /2$ is odd, so that $\left(
-1\right)  ^{n\left(  n-1\right)  /2}=-1$. Hence,%
\[
\left(  -1\right)  ^{w_{0}}=\left(  -1\right)  ^{n\left(  n-1\right)  /2}=-1.
\]
Compared with $%
\begin{cases}
1, & \text{if }n\equiv0\operatorname{mod}4\text{ or }n\equiv
1\operatorname{mod}4;\\
-1, & \text{if }n\equiv2\operatorname{mod}4\text{ or }n\equiv
3\operatorname{mod}4
\end{cases}
=-1$ (since $n\equiv2\operatorname{mod}4$ or $n\equiv3\operatorname{mod}4$
(since $n\equiv3\operatorname{mod}4$)), this yields
\[
\left(  -1\right)  ^{w_{0}}=
\begin{cases}
1, & \text{if }n\equiv0\operatorname{mod}4\text{ or }n\equiv
1\operatorname{mod}4;\\
-1, & \text{if }n\equiv2\operatorname{mod}4\text{ or }n\equiv
3\operatorname{mod}4
\end{cases}
.
\]
Thus, (\ref{sol.ps4.1c.extra}) is proven in Case 4.

We now have proven (\ref{sol.ps4.1c.extra}) in each of the four Cases 1, 2, 3
and 4. Thus, (\ref{sol.ps4.1c.extra}) always holds.]

This finishes our solution of Exercise \ref{exe.ps4.1c}.
\end{proof}
\end{verlong}

\subsection{Solution to Exercise \ref{exe.ps4.2}}

\begin{proof}
[Solution to Exercise \ref{exe.ps4.2}.]\textbf{(a)} Let $\sigma$ be a
permutation of $X$. We need to prove that $\left(  -1\right)  _{\phi}^{\sigma
}$ depends only on the permutation $\sigma$ of $X$, but not on the bijection
$\phi$. In other words, we need to prove that any two different choices of
$\phi$ will lead to the same $\left(  -1\right)  _{\phi}^{\sigma}$. In other
words, we need to prove that if $\phi_{1}$ and $\phi_{2}$ are two bijections
$\phi:X\rightarrow\left\{  1,2,\ldots,n\right\}  $ for some $n\in\mathbb{N}$
(possibly distinct), then $\left(  -1\right)  _{\phi_{1}}^{\sigma}=\left(
-1\right)  _{\phi_{2}}^{\sigma}$.

So let $\phi_{1}$ and $\phi_{2}$ be two bijections $\phi:X\rightarrow\left\{
1,2,\ldots,n\right\}  $ for some $n\in\mathbb{N}$ (possibly distinct). We must
show that $\left(  -1\right)  _{\phi_{1}}^{\sigma}=\left(  -1\right)
_{\phi_{2}}^{\sigma}$.

The map $\sigma$ is a permutation of $X$, thus a bijection $X\rightarrow X$.

\begin{vershort}
We know that $\phi_{1}$ is a bijection $X\rightarrow\left\{  1,2,\ldots
,n\right\}  $ for some $n\in\mathbb{N}$. In other words, there exists some
$n\in\mathbb{N}$ such that $\phi_{1}$ is a bijection $X\rightarrow\left\{
1,2,\ldots,n\right\}  $. Denote this $n$ by $n_{1}$. Thus, $\phi_{1}$ is a
bijection $X\rightarrow\left\{  1,2,\ldots,n_{1}\right\}  $. The definition of
$\left(  -1\right)  _{\phi_{1}}^{\sigma}$ yields $\left(  -1\right)
_{\phi_{1}}^{\sigma}=\left(  -1\right)  ^{\phi_{1}\circ\sigma\circ\phi
_{1}^{-1}}$.
\end{vershort}

\begin{verlong}
We know that $\phi_{1}$ is a bijection $\phi:X\rightarrow\left\{
1,2,\ldots,n\right\}  $ for some $n\in\mathbb{N}$. In other words, $\phi_{1}$
is a bijection $X\rightarrow\left\{  1,2,\ldots,n\right\}  $ for some
$n\in\mathbb{N}$. In other words, there exists some $n\in\mathbb{N}$ such that
$\phi_{1}$ is a bijection $X\rightarrow\left\{  1,2,\ldots,n\right\}  $.
Denote this $n$ by $n_{1}$. Thus, $\phi_{1}$ is a bijection $X\rightarrow
\left\{  1,2,\ldots,n_{1}\right\}  $. The definition of $\left(  -1\right)
_{\phi_{1}}^{\sigma}$ yields $\left(  -1\right)  _{\phi_{1}}^{\sigma}=\left(
-1\right)  ^{\phi_{1}\circ\sigma\circ\phi_{1}^{-1}}$.
\end{verlong}

The map $\phi_{1}$ is a bijection. Thus, its inverse $\phi_{1}^{-1}$ is
well-defined and also a bijection.

The map $\phi_{1}\circ\sigma\circ\phi_{1}^{-1}:\left\{  1,2,\ldots
,n_{1}\right\}  \rightarrow\left\{  1,2,\ldots,n_{1}\right\}  $ is a bijection
(since it is a composition of the three bijections $\phi_{1}$, $\sigma$ and
$\phi_{1}^{-1}$). In other words, the map $\phi_{1}\circ\sigma\circ\phi
_{1}^{-1}$ is a permutation of $\left\{  1,2,\ldots,n_{1}\right\}  $. Thus,
$\phi_{1}\circ\sigma\circ\phi_{1}^{-1}\in S_{n_{1}}$.

\begin{vershort}
We thus have shown that $\phi_{1}^{-1}$ is well-defined and a bijection, and
found an $n_{1}\in\mathbb{N}$ such that $\phi_{1}\circ\sigma\circ\phi_{1}%
^{-1}\in S_{n_{1}}$. Similarly, we can show that $\phi_{2}^{-1}$ is
well-defined and a bijection, and find an $n_{2}\in\mathbb{N}$ such that
$\phi_{2}\circ\sigma\circ\phi_{2}^{-1}\in S_{n_{2}}$. Consider this $n_{2}$.
(We shall soon see that $n_{1}=n_{2}$.) We have $\left(  -1\right)  _{\phi
_{2}}^{\sigma}=\left(  -1\right)  ^{\phi_{2}\circ\sigma\circ\phi_{2}^{-1}}$
(by the definition of $\left(  -1\right)  _{\phi_{2}}^{\sigma}$).
\end{vershort}

\begin{verlong}
We know that $\phi_{2}$ is a bijection $\phi:X\rightarrow\left\{
1,2,\ldots,n\right\}  $ for some $n\in\mathbb{N}$. In other words, $\phi_{2}$
is a bijection $X\rightarrow\left\{  1,2,\ldots,n\right\}  $ for some
$n\in\mathbb{N}$. In other words, there exists some $n\in\mathbb{N}$ such that
$\phi_{2}$ is a bijection $X\rightarrow\left\{  1,2,\ldots,n\right\}  $.
Denote this $n$ by $n_{2}$. Thus, $\phi_{2}$ is a bijection $X\rightarrow
\left\{  1,2,\ldots,n_{2}\right\}  $. The definition of $\left(  -1\right)
_{\phi_{2}}^{\sigma}$ yields $\left(  -1\right)  _{\phi_{2}}^{\sigma}=\left(
-1\right)  ^{\phi_{2}\circ\sigma\circ\phi_{2}^{-1}}$. (We will soon see that
$n_{1}=n_{2}$.)

The map $\phi_{2}$ is a bijection. Thus, its inverse $\phi_{2}^{-1}$ is
well-defined and also a bijection.

The map $\phi_{2}\circ\sigma\circ\phi_{2}^{-1}:\left\{  1,2,\ldots
,n_{2}\right\}  \rightarrow\left\{  1,2,\ldots,n_{2}\right\}  $ is a bijection
(since it is a composition of the three bijections $\phi_{2}$, $\sigma$ and
$\phi_{2}^{-1}$). In other words, the map $\phi_{2}\circ\sigma\circ\phi
_{2}^{-1}$ is a permutation of $\left\{  1,2,\ldots,n_{2}\right\}  $. Thus,
$\phi_{2}\circ\sigma\circ\phi_{2}^{-1}\in S_{n_{2}}$.
\end{verlong}

There exists a bijection from $X$ to $\left\{  1,2,\ldots,n_{1}\right\}  $
(namely, $\phi_{1}$). Hence, $\left\vert X\right\vert =\left\vert \left\{
1,2,\ldots,n_{1}\right\}  \right\vert =n_{1}$. Similarly, $\left\vert
X\right\vert =n_{2}$. Thus, $n_{1}=\left\vert X\right\vert =n_{2}$. Thus, we
can define an $n\in\mathbb{N}$ by $n=n_{1}=n_{2}$. Consider this $n$.

The map $\phi_{2}\circ\phi_{1}^{-1}:\left\{  1,2,\ldots,n_{1}\right\}
\rightarrow\left\{  1,2,\ldots,n_{2}\right\}  $ is a bijection (since it is a
composition of two bijections). Since $n_{1}=n$ and $n_{2}=n$, this rewrites
as follows: The map $\phi_{2}\circ\phi_{1}^{-1}:\left\{  1,2,\ldots,n\right\}
\rightarrow\left\{  1,2,\ldots,n\right\}  $ is a bijection. In other words,
the map $\phi_{2}\circ\phi_{1}^{-1}$ is a permutation of $\left\{
1,2,\ldots,n\right\}  $. Thus, $\phi_{2}\circ\phi_{1}^{-1}\in S_{n}$.

We have $\phi_{1}\circ\sigma\circ\phi_{1}^{-1}\in S_{n_{1}}=S_{n}$ (since
$n_{1}=n$) and $\phi_{2}\circ\sigma\circ\phi_{2}^{-1}\in S_{n_{2}}=S_{n}$
(since $n_{2}=n$).

Now, (\ref{eq.sign.prod}) (applied to $\phi_{2}\circ\phi_{1}^{-1}$ and
$\phi_{1}\circ\sigma\circ\phi_{1}^{-1}$ instead of $\sigma$ and $\tau$) yields%
\[
\left(  -1\right)  ^{\phi_{2}\circ\phi_{1}^{-1}\circ\phi_{1}\circ\sigma
\circ\phi_{1}^{-1}}=\left(  -1\right)  ^{\phi_{2}\circ\phi_{1}^{-1}}%
\cdot\left(  -1\right)  ^{\phi_{1}\circ\sigma\circ\phi_{1}^{-1}}.
\]
Since $\phi_{2}\circ\underbrace{\phi_{1}^{-1}\circ\phi_{1}}%
_{=\operatorname*{id}}\circ\sigma\circ\phi_{1}^{-1}=\phi_{2}\circ\sigma
\circ\phi_{1}^{-1}$, this rewrites as
\begin{equation}
\left(  -1\right)  ^{\phi_{2}\circ\sigma\circ\phi_{1}^{-1}}=\left(  -1\right)
^{\phi_{2}\circ\phi_{1}^{-1}}\cdot\left(  -1\right)  ^{\phi_{1}\circ
\sigma\circ\phi_{1}^{-1}}. \label{sol.ps4.2.1}%
\end{equation}

On the other hand, (\ref{eq.sign.prod}) (applied to $\phi_{2}\circ\sigma
\circ\phi_{2}^{-1}$ and $\phi_{2}\circ\phi_{1}^{-1}$ instead of $\sigma$ and
$\tau$) yields%
\[
\left(  -1\right)  ^{\phi_{2}\circ\sigma\circ\phi_{2}^{-1}\circ\phi_{2}%
\circ\phi_{1}^{-1}}=\left(  -1\right)  ^{\phi_{2}\circ\sigma\circ\phi_{2}%
^{-1}}\cdot\left(  -1\right)  ^{\phi_{2}\circ\phi_{1}^{-1}}=\left(  -1\right)
^{\phi_{2}\circ\phi_{1}^{-1}}\cdot\left(  -1\right)  ^{\phi_{2}\circ
\sigma\circ\phi_{2}^{-1}}.
\]
Since $\phi_{2}\circ\sigma\circ\underbrace{\phi_{2}^{-1}\circ\phi_{2}%
}_{=\operatorname*{id}}\circ\phi_{1}^{-1}=\phi_{2}\circ\sigma\circ\phi
_{1}^{-1}$, this rewrites as
\[
\left(  -1\right)  ^{\phi_{2}\circ\sigma\circ\phi_{1}^{-1}}=\left(  -1\right)
^{\phi_{2}\circ\phi_{1}^{-1}}\cdot\left(  -1\right)  ^{\phi_{2}\circ
\sigma\circ\phi_{2}^{-1}}.
\]
Comparing this with (\ref{sol.ps4.2.1}), we obtain
\[
\left(  -1\right)  ^{\phi_{2}\circ\phi_{1}^{-1}}\cdot\left(  -1\right)
^{\phi_{1}\circ\sigma\circ\phi_{1}^{-1}}=\left(  -1\right)  ^{\phi_{2}%
\circ\phi_{1}^{-1}}\cdot\left(  -1\right)  ^{\phi_{2}\circ\sigma\circ\phi
_{2}^{-1}}.
\]
We can cancel $\left(  -1\right)  ^{\phi_{2}\circ\phi_{1}^{-1}}$ from this
equality (since $\left(  -1\right)  ^{\phi_{2}\circ\phi_{1}^{-1}}\in\left\{
1,-1\right\}  $ is a nonzero integer), and thus obtain $\left(  -1\right)
^{\phi_{1}\circ\sigma\circ\phi_{1}^{-1}}=\left(  -1\right)  ^{\phi_{2}%
\circ\sigma\circ\phi_{2}^{-1}}$. Hence,%
\[
\left(  -1\right)  _{\phi_{1}}^{\sigma}=\left(  -1\right)  ^{\phi_{1}%
\circ\sigma\circ\phi_{1}^{-1}}=\left(  -1\right)  ^{\phi_{2}\circ\sigma
\circ\phi_{2}^{-1}}=\left(  -1\right)  _{\phi_{2}}^{\sigma}.
\]
As we know, this completes the solution of Exercise \ref{exe.ps4.2}
\textbf{(a)}.

\textbf{(b)} Fix a bijection $\phi:X\rightarrow\left\{  1,2,\ldots,n\right\}
$ for some $n\in\mathbb{N}$. (Such a bijection always exists.) We shall denote
the identity permutation of $X$ by $\operatorname*{id}\nolimits_{X}$, so as to
distinguish it from the identity permutation of $\left\{  1,2,\ldots
,n\right\}  $ (which we keep denoting by $\operatorname*{id}$ without a
subscript). The definition of $\left(  -1\right)  ^{\operatorname*{id}%
\nolimits_{X}}$ now yields
\begin{align*}
\left(  -1\right)  ^{\operatorname*{id}\nolimits_{X}}  &  =\left(  -1\right)
_{\phi}^{\operatorname*{id}\nolimits_{X}}=\left(  -1\right)  ^{\phi
\circ\operatorname*{id}\nolimits_{X}\circ\phi^{-1}}\ \ \ \ \ \ \ \ \ \ \left(
\text{by the definition of }\left(  -1\right)  _{\phi}^{\operatorname*{id}%
\nolimits_{X}}\right) \\
&  =\left(  -1\right)  ^{\operatorname*{id}}\ \ \ \ \ \ \ \ \ \ \left(
\text{since }\phi\circ\operatorname*{id}\nolimits_{X}\circ\phi^{-1}=\phi
\circ\phi^{-1}=\operatorname*{id}\right) \\
&  =1.
\end{align*}
In other words, $\left(  -1\right)  ^{\operatorname*{id}}=1$ for the identity
permutation $\operatorname*{id}:X\rightarrow X$ of $X$. Exercise
\ref{exe.ps4.2} \textbf{(b)} is thus solved.

\textbf{(c)} Let $\sigma$ and $\tau$ be two permutations of $X$. Fix a
bijection $\phi:X\rightarrow\left\{  1,2,\ldots,n\right\}  $ for some
$n\in\mathbb{N}$. (Such a bijection always exists.) The definition of $\left(
-1\right)  ^{\sigma}$ yields%
\begin{equation}
\left(  -1\right)  ^{\sigma}=\left(  -1\right)  _{\phi}^{\sigma}=\left(
-1\right)  ^{\phi\circ\sigma\circ\phi^{-1}}\ \ \ \ \ \ \ \ \ \ \left(
\text{by the definition of }\left(  -1\right)  _{\phi}^{\sigma}\right)  .
\label{sol.ps4.2.c.1}%
\end{equation}
The definition of $\left(  -1\right)  ^{\tau}$ yields%
\begin{equation}
\left(  -1\right)  ^{\tau}=\left(  -1\right)  _{\phi}^{\tau}=\left(
-1\right)  ^{\phi\circ\tau\circ\phi^{-1}}\ \ \ \ \ \ \ \ \ \ \left(  \text{by
the definition of }\left(  -1\right)  _{\phi}^{\tau}\right)  .
\label{sol.ps4.2.c.2}%
\end{equation}
The maps $\phi\circ\sigma\circ\phi^{-1}$ and $\phi\circ\tau\circ\phi^{-1}$ are
permutations in $S_{n}$. Therefore, (\ref{eq.sign.prod}) (applied to
$\phi\circ\sigma\circ\phi^{-1}$ and $\phi\circ\tau\circ\phi^{-1}$ instead of
$\sigma$ and $\tau$) yields%
\[
\left(  -1\right)  ^{\phi\circ\sigma\circ\phi^{-1}\circ\phi\circ\tau\circ
\phi^{-1}}=\left(  -1\right)  ^{\phi\circ\sigma\circ\phi^{-1}}\cdot\left(
-1\right)  ^{\phi\circ\tau\circ\phi^{-1}}.
\]
Since $\phi\circ\sigma\circ\underbrace{\phi^{-1}\circ\phi}%
_{=\operatorname*{id}}\circ\tau\circ\phi^{-1}=\phi\circ\sigma\circ\tau
\circ\phi^{-1}$, this rewrites as
\[
\left(  -1\right)  ^{\phi\circ\sigma\circ\tau\circ\phi^{-1}}=\left(
-1\right)  ^{\phi\circ\sigma\circ\phi^{-1}}\cdot\left(  -1\right)  ^{\phi
\circ\tau\circ\phi^{-1}}.
\]
But the definition of $\left(  -1\right)  ^{\sigma\circ\tau}$ yields%
\begin{align*}
\left(  -1\right)  ^{\sigma\circ\tau}  &  =\left(  -1\right)  _{\phi}%
^{\sigma\circ\tau}=\left(  -1\right)  ^{\phi\circ\sigma\circ\tau\circ\phi
^{-1}}\ \ \ \ \ \ \ \ \ \ \left(  \text{by the definition of }\left(
-1\right)  _{\phi}^{\sigma\circ\tau}\right) \\
&  =\underbrace{\left(  -1\right)  ^{\phi\circ\sigma\circ\phi^{-1}}%
}_{\substack{=\left(  -1\right)  ^{\sigma}\\\text{(by (\ref{sol.ps4.2.c.1}))}%
}}\cdot\underbrace{\left(  -1\right)  ^{\phi\circ\tau\circ\phi^{-1}}%
}_{\substack{=\left(  -1\right)  ^{\tau}\\\text{(by (\ref{sol.ps4.2.c.2}))}%
}}=\left(  -1\right)  ^{\sigma}\cdot\left(  -1\right)  ^{\tau}.
\end{align*}
This solves Exercise \ref{exe.ps4.2} \textbf{(c)}.
\end{proof}

\subsection{Solution to Exercise \ref{exe.perm.sign.pseudoexplicit}}

\begin{vershort}
\begin{proof}
[Solution to Exercise \ref{exe.perm.sign.pseudoexplicit}.]\textbf{(b)} We
start with some trivia on sets and subsets.

If $A$ is any set, then $\mathcal{P}_{2}\left(  A\right)  $ shall denote the
set of all $2$-element subsets of $A$. In other words, $\mathcal{P}_{2}\left(
A\right)  $ is defined to be $\left\{  S\subseteq A\ \mid\ \left\vert
S\right\vert =2\right\}  $. For instance,
\begin{align*}
\mathcal{P}_{2}\left(  \left\{  3,6,7\right\}  \right)   &  =\left\{  \left\{
3,6\right\}  ,\left\{  3,7\right\}  ,\left\{  6,7\right\}  \right\}  ;\\
\mathcal{P}_{2}\left(  \left\{  1,2,3,4\right\}  \right)   &  =\left\{
\left\{  1,2\right\}  ,\left\{  1,3\right\}  ,\left\{  1,4\right\}  ,\left\{
2,3\right\}  ,\left\{  2,4\right\}  ,\left\{  3,4\right\}  \right\}  ;\\
\mathcal{P}_{2}\left(  \left\{  3\right\}  \right)   &  =\varnothing;\\
\mathcal{P}_{2}\left(  \varnothing\right)   &  =\varnothing.
\end{align*}
\footnote{Several authors write $\dbinom{A}{2}$ for $\mathcal{P}_{2}\left(
A\right)  $. This notation looks like a binomial coefficient, but with $A$ at
the top. Of course, this notation is chosen for its suggestiveness: When $A$
is finite, the set $\dbinom{A}{2}$ satisfies $\left\vert \dbinom{A}%
{2}\right\vert =\dbinom{\left\vert A\right\vert }{2}$.}

If $A$ and $B$ are two sets, and if $f:A\rightarrow B$ is an injective map,
then we can define a map $f_{\ast}:\mathcal{P}_{2}\left(  A\right)
\rightarrow\mathcal{P}_{2}\left(  B\right)  $ by%
\[
\left(  f_{\ast}\left(  S\right)  =f\left(  S\right)
\ \ \ \ \ \ \ \ \ \ \text{for every }S\in\mathcal{P}_{2}\left(  A\right)
\right)  .
\]
\footnote{At this point, we need to check that this map $f_{\ast}$ is
well-defined. Before I do this, let me rewrite the definition of $f_{\ast}$ in
a more intuitive way: An element of $\mathcal{P}_{2}\left(  A\right)  $ is a
$2$-element subset $\left\{  a,a^{\prime}\right\}  $ of $A$. The map $f_{\ast
}$ takes this subset to $\left\{  f\left(  a\right)  ,f\left(  a^{\prime
}\right)  \right\}  $ (in other words, it applies $f$ to each of its
elements).
\par
So why is the map $f_{\ast}$ well-defined? It is supposed to send every
$S\in\mathcal{P}_{2}\left(  A\right)  $ to $f\left(  S\right)  $. Thus, in
order to prove that it is well-defined, we need to show that $f\left(
S\right)  \in\mathcal{P}_{2}\left(  B\right)  $ for every $S\in\mathcal{P}%
_{2}\left(  A\right)  $.
\par
Let $S\in\mathcal{P}_{2}\left(  A\right)  $. Thus, the set $S$ is a
$2$-element subset of $A$. The map $f$ sends its two elements to two
\textbf{distinct} elements of $B$ (they are distinct because $f$ is
injective). In other words, the set $f\left(  S\right)  $ has $2$ elements.
Thus, $f\left(  S\right)  $ is a $2$-element subset of $B$; in other words,
$f\left(  S\right)  \in\mathcal{P}_{2}\left(  B\right)  $. This proves that
the map $f_{\ast}$ is well-defined.
\par
Notice that we have used the injectivity of $f$ in this argument.} This
construction has the following three basic properties:

\begin{enumerate}
\item If $A$ is a set, then $\left(  \operatorname*{id}\nolimits_{A}\right)
_{\ast}=\operatorname*{id}\nolimits_{\mathcal{P}_{2}\left(  A\right)  }$.

\item If $A$, $B$ and $C$ are three sets, and if $f:A\rightarrow B$ and
$g:B\rightarrow C$ are two injective maps, then $\left(  g\circ f\right)
_{\ast}=g_{\ast}\circ f_{\ast}$. (Of course, $g\circ f$ is injective here, so
$\left(  g\circ f\right)  _{\ast}$ makes sense.)

\item If $A$ and $B$ are two sets, and if $f:A\rightarrow B$ is an invertible
map, then $f_{\ast}$ is invertible as well and satisfies $\left(
f^{-1}\right)  _{\ast}=\left(  f_{\ast}\right)  ^{-1}$.
\end{enumerate}

(The first of these three properties is obvious; the second follows by
observing that $\left(  g\circ f\right)  \left(  S\right)  =g\left(  f\left(
S\right)  \right)  $ for every $S\in\mathcal{P}_{2}\left(  A\right)  $; the
third can be proven using the second or directly.)

Let $\left[  n\right]  $ be the set $\left\{  1,2,\ldots,n\right\}  $. Recall
that $S_{n}$ is the set of all permutations of the set $\left\{
1,2,\ldots,n\right\}  $. In other words, $S_{n}$ is the set of all
permutations of the set $\left[  n\right]  $ (since $\left\{  1,2,\ldots
,n\right\}  =\left[  n\right]  $).

We have $\sigma\in S_{n}$. In other words, $\sigma$ is a permutation of the
set $\left[  n\right]  $ (since $S_{n}$ is the set of all such permutations).
Thus, $\sigma$ is a bijective map $\left[  n\right]  \rightarrow\left[
n\right]  $. In particular, the map $\sigma_{\ast}:\mathcal{P}_{2}\left(
\left[  n\right]  \right)  \rightarrow\mathcal{P}_{2}\left(  \left[  n\right]
\right)  $ is well-defined.

Recall that if $A$ and $B$ are two sets, and if $f:A\rightarrow B$ is an
invertible map, then $f_{\ast}$ is invertible as well. Applying this to
$A=\left[  n\right]  $, $B=\left[  n\right]  $ and $f=\sigma$, we conclude
that $\sigma_{\ast}$ is invertible (since $\sigma$ is invertible). In other
words, $\sigma_{\ast}$ is bijective.

Let $G$ be the subset%
\[
\left\{  \left(  u,v\right)  \in\left[  n\right]  ^{2}\ \mid\ u<v\right\}
\]
of $\left[  n\right]  ^{2}$.

For example, if $n=4$, then%
\[
G=\left\{  \left(  1,2\right)  ,\left(  1,3\right)  ,\left(  1,4\right)
,\left(  2,3\right)  ,\left(  2,4\right)  ,\left(  3,4\right)  \right\}  .
\]
Comparing this with%
\[
\mathcal{P}_{2}\left(  \left[  4\right]  \right)  =\left\{  \left\{
1,2\right\}  ,\left\{  1,3\right\}  ,\left\{  1,4\right\}  ,\left\{
2,3\right\}  ,\left\{  2,4\right\}  ,\left\{  3,4\right\}  \right\}  ,
\]
we observe that the set $\mathcal{P}_{2}\left(  \left[  n\right]  \right)  $
is obtained from $G$ by \textquotedblleft replacing all parentheses by
brackets\textquotedblright\ (i.e., replacing each $\left(  i,j\right)  \in G$
by $\left\{  i,j\right\}  $). This holds for all $n$, not just for $n=4$. Let
us make this observation somewhat more bulletproof: Each $\left(  i,j\right)
\in G$ satisfies $\left\{  i,j\right\}  \in\mathcal{P}_{2}\left(  \left[
n\right]  \right)  $ (because $\left(  i,j\right)  \in G$ leads to $i<j$, so
that $i\neq j$, so that $\left\{  i,j\right\}  $ is a $2$-element set, and
thus $\left\{  i,j\right\}  \in\mathcal{P}_{2}\left(  \left[  n\right]
\right)  $). Thus, we can define a map $\rho:G\rightarrow\mathcal{P}%
_{2}\left(  \left[  n\right]  \right)  $ by setting%
\[
\left(  \rho\left(  \left(  i,j\right)  \right)  =\left\{  i,j\right\}
\ \ \ \ \ \ \ \ \ \ \text{for every }\left(  i,j\right)  \in G\right)  .
\]
This map $\rho$ is injective (indeed, we can reconstruct every $\left(
i,j\right)  \in G$ from its image $\rho\left(  \left(  i,j\right)  \right)
=\left\{  i,j\right\}  $, because $\left(  i,j\right)  \in G$ entails $i<j$)
and surjective (since every two-element subset $S$ of $\left[  n\right]  $ has
the form $\left\{  i,j\right\}  $ for some $\left(  i,j\right)  \in\left[
n\right]  ^{2}$ satisfying $i<j$). Hence, the map $\rho$ is bijective.

For every $S\in\mathcal{P}_{2}\left(  \left[  n\right]  \right)  $, we define
an element $a_{S}$ of $\mathbb{C}$ by%
\[
a_{S}=a_{\left(  \min S,\max S\right)  }.
\]
(In other words, for every $S\in\mathcal{P}_{2}\left(  \left[  n\right]
\right)  $, we define an element $a_{S}$ of $\mathbb{C}$ by $a_{S}=a_{\left(
i,j\right)  }$, where $i$ and $j$ are the two elements of $S$ in increasing order.)

Let $\operatorname*{Inv}\left(  \sigma\right)  $ be the set of inversions of
$\sigma$. Then, $\left(  -1\right)  ^{\sigma}=\left(  -1\right)  ^{\left\vert
\operatorname*{Inv}\left(  \sigma\right)  \right\vert }$%
\ \ \ \ \footnote{\textit{Proof.} The definition of $\ell\left(
\sigma\right)  $ shows that%
\begin{align*}
\ell\left(  \sigma\right)   &  =\left(  \text{the number of inversions of
}\sigma\right)  =\left(  \text{the number of elements of }\operatorname*{Inv}%
\left(  \sigma\right)  \right) \\
&  \ \ \ \ \ \ \ \ \ \ \left(  \text{since }\operatorname*{Inv}\left(
\sigma\right)  \text{ is the set of all inversions of }\sigma\right) \\
&  =\left\vert \operatorname*{Inv}\left(  \sigma\right)  \right\vert .
\end{align*}
\par
But the definition of $\left(  -1\right)  ^{\sigma}$ yields $\left(
-1\right)  ^{\sigma}=\left(  -1\right)  ^{\ell\left(  \sigma\right)  }=\left(
-1\right)  ^{\left\vert \operatorname*{Inv}\left(  \sigma\right)  \right\vert
}$ (since $\ell\left(  \sigma\right)  =\left\vert \operatorname*{Inv}\left(
\sigma\right)  \right\vert $), qed.}. Furthermore, $\operatorname*{Inv}\left(
\sigma\right)  \subseteq G$\ \ \ \ \footnote{\textit{Proof.} Let
$c\in\operatorname*{Inv}\left(  \sigma\right)  $. Thus, $c$ is an inversion of
$\sigma$ (since $\operatorname*{Inv}\left(  \sigma\right)  $ is the set of
inversions of $\sigma$). In other words, $c$ is a pair $\left(  i,j\right)  $
of integers satisfying $1\leq i<j\leq n$ and $\sigma\left(  i\right)
>\sigma\left(  j\right)  $. Consider this $\left(  i,j\right)  $. We have
$i\in\left[  n\right]  $ (since $1\leq i\leq n$) and $j\in\left[  n\right]  $
(since $1\leq j\leq n$). Thus, $\left(  i,j\right)  \in\left[  n\right]  ^{2}%
$. Thus, $\left(  i,j\right)  $ is an element of $\left[  n\right]  ^{2}$ and
satisfies $i<j$. In other words, $\left(  i,j\right)  \in\left\{  \left(
u,v\right)  \in\left[  n\right]  ^{2}\ \mid\ u<v\right\}  =G$. Hence,
$c=\left(  i,j\right)  \in G$.
\par
Now, let us forget that we fixed $c$. We thus have proven that every
$c\in\operatorname*{Inv}\left(  \sigma\right)  $ satisfies $c\in G$. In other
words, $\operatorname*{Inv}\left(  \sigma\right)  \subseteq G$, qed.}.

Now, we notice the following facts:

\begin{itemize}
\item For every $\left(  i,j\right)  \in G$, we have%
\begin{equation}
a_{\left(  i,j\right)  }=a_{\rho\left(  \left(  i,j\right)  \right)  }.
\label{sol.perm.sign.pseudoexplicit.short.b.a1}%
\end{equation}
\footnote{\textit{Proof of (\ref{sol.perm.sign.pseudoexplicit.short.b.a1}):}
Let $\left(  i,j\right)  \in G$. Thus, $\left(  i,j\right)  \in G=\left\{
\left(  u,v\right)  \in\left[  n\right]  ^{2}\ \mid\ u<v\right\}  $. In other
words, $\left(  i,j\right)  \in\left[  n\right]  ^{2}$ and $i<j$.
\par
The definition of $\rho$ shows that $\rho\left(  \left(  i,j\right)  \right)
=\left\{  i,j\right\}  $. Since $i<j$, this shows that the elements of the set
$\rho\left(  \left(  i,j\right)  \right)  $ listed in increasing order are $i$
and $j$. Hence, $\min\left(  \rho\left(  \left(  i,j\right)  \right)  \right)
=i$ and $\max\left(  \rho\left(  \left(  i,j\right)  \right)  \right)  =j$.
\par
Now, the definition of $a_{\rho\left(  \left(  i,j\right)  \right)  }$ shows
that $a_{\rho\left(  \left(  i,j\right)  \right)  }=a_{\left(  \min\left(
\rho\left(  \left(  i,j\right)  \right)  \right)  ,\max\left(  \rho\left(
\left(  i,j\right)  \right)  \right)  \right)  }=a_{\left(  i,j\right)  }$
(since $\min\left(  \rho\left(  \left(  i,j\right)  \right)  \right)  =i$ and
$\max\left(  \rho\left(  \left(  i,j\right)  \right)  \right)  =j$). This
proves (\ref{sol.perm.sign.pseudoexplicit.short.b.a1}).}

\item For every $\left(  i,j\right)  \in\operatorname*{Inv}\left(
\sigma\right)  $, we have%
\begin{equation}
a_{\left(  \sigma\left(  i\right)  ,\sigma\left(  j\right)  \right)
}=-a_{\sigma\left(  \rho\left(  \left(  i,j\right)  \right)  \right)  }.
\label{sol.perm.sign.pseudoexplicit.short.b.a2}%
\end{equation}
\footnote{\textit{Proof of (\ref{sol.perm.sign.pseudoexplicit.short.b.a2}):}
Let $\left(  i,j\right)  \in\operatorname*{Inv}\left(  \sigma\right)  $. Thus,
$\left(  i,j\right)  $ is an inversion of $\sigma$. In other words, $\left(
i,j\right)  $ is a pair of integers satisfying $1\leq i<j\leq n$ and
$\sigma\left(  i\right)  >\sigma\left(  j\right)  $.
\par
The definition of $\rho$ shows that $\rho\left(  \left(  i,j\right)  \right)
=\left\{  i,j\right\}  $. Hence, $\sigma\left(  \underbrace{\rho\left(
\left(  i,j\right)  \right)  }_{=\left\{  i,j\right\}  }\right)
=\sigma\left(  \left\{  i,j\right\}  \right)  =\left\{  \sigma\left(
i\right)  ,\sigma\left(  j\right)  \right\}  $. Since $\sigma\left(  i\right)
>\sigma\left(  j\right)  $, this shows that the elements of the set
$\sigma\left(  \rho\left(  \left(  i,j\right)  \right)  \right)  $ listed in
increasing order are $\sigma\left(  j\right)  $ and $\sigma\left(  i\right)
$. Hence, $\min\left(  \sigma\left(  \rho\left(  \left(  i,j\right)  \right)
\right)  \right)  =\sigma\left(  j\right)  $ and $\max\left(  \sigma\left(
\rho\left(  \left(  i,j\right)  \right)  \right)  \right)  =\sigma\left(
i\right)  $.
\par
Now, the definition of $a_{\sigma\left(  \rho\left(  \left(  i,j\right)
\right)  \right)  }$ shows that%
\begin{align*}
a_{\sigma\left(  \rho\left(  \left(  i,j\right)  \right)  \right)  }  &
=a_{\left(  \min\left(  \sigma\left(  \rho\left(  \left(  i,j\right)  \right)
\right)  \right)  ,\max\left(  \sigma\left(  \rho\left(  \left(  i,j\right)
\right)  \right)  \right)  \right)  }=a_{\left(  \sigma\left(  j\right)
,\sigma\left(  i\right)  \right)  }\\
&  \ \ \ \ \ \ \ \ \ \ \left(  \text{since }\min\left(  \sigma\left(
\rho\left(  \left(  i,j\right)  \right)  \right)  \right)  =\sigma\left(
j\right)  \text{ and }\max\left(  \sigma\left(  \rho\left(  \left(
i,j\right)  \right)  \right)  \right)  =\sigma\left(  i\right)  \right) \\
&  =-a_{\left(  \sigma\left(  i\right)  ,\sigma\left(  j\right)  \right)
}\ \ \ \ \ \ \ \ \ \ \left(  \text{by
(\ref{eq.exe.perm.sign.pseudoexplicit.b.skew}), applied to }\sigma\left(
i\right)  \text{ and }\sigma\left(  j\right)  \text{ instead of }i\text{ and
}j\right)  .
\end{align*}
Hence, $a_{\left(  \sigma\left(  i\right)  ,\sigma\left(  j\right)  \right)
}=-a_{\sigma\left(  \rho\left(  \left(  i,j\right)  \right)  \right)  }$. This
proves (\ref{sol.perm.sign.pseudoexplicit.short.b.a2}).}

\item For every $\left(  i,j\right)  \in G$ satisfying $\left(  i,j\right)
\notin\operatorname*{Inv}\left(  \sigma\right)  $, we have%
\begin{equation}
a_{\left(  \sigma\left(  i\right)  ,\sigma\left(  j\right)  \right)
}=a_{\sigma\left(  \rho\left(  \left(  i,j\right)  \right)  \right)  }.
\label{sol.perm.sign.pseudoexplicit.short.b.a3}%
\end{equation}
\footnote{\textit{Proof of (\ref{sol.perm.sign.pseudoexplicit.short.b.a3}):}
Let $\left(  i,j\right)  \in G$ be such that $\left(  i,j\right)
\notin\operatorname*{Inv}\left(  \sigma\right)  $. Thus, $\left(  i,j\right)
$ is an element of $\left[  n\right]  ^{2}$ satisfying $i<j$ (since $\left(
i,j\right)  \in G=\left\{  \left(  u,v\right)  \in\left[  n\right]  ^{2}%
\ \mid\ u<v\right\}  $).
\par
If we had $\sigma\left(  i\right)  >\sigma\left(  j\right)  $, then $\left(
i,j\right)  $ would be an inversion of $\sigma$ (since $\left(  i,j\right)  $
is a pair of integers satisfying $1\leq i<j\leq n$), and thus would belong to
$\operatorname*{Inv}\left(  \sigma\right)  $; this would contradict $\left(
i,j\right)  \notin\operatorname*{Inv}\left(  \sigma\right)  $. Hence, we
cannot have $\sigma\left(  i\right)  >\sigma\left(  j\right)  $. Thus, we have
$\sigma\left(  i\right)  \leq\sigma\left(  j\right)  $. Since $\sigma\left(
i\right)  \neq\sigma\left(  j\right)  $ (because $i\neq j$ and because
$\sigma$ is injective), this shows that $\sigma\left(  i\right)
<\sigma\left(  j\right)  $.
\par
The definition of $\rho$ shows that $\rho\left(  \left(  i,j\right)  \right)
=\left\{  i,j\right\}  $. Hence, $\sigma\left(  \underbrace{\rho\left(
\left(  i,j\right)  \right)  }_{=\left\{  i,j\right\}  }\right)
=\sigma\left(  \left\{  i,j\right\}  \right)  =\left\{  \sigma\left(
i\right)  ,\sigma\left(  j\right)  \right\}  $. Since $\sigma\left(  i\right)
<\sigma\left(  j\right)  $, this shows that the elements of the set
$\sigma\left(  \rho\left(  \left(  i,j\right)  \right)  \right)  $ listed in
increasing order are $\sigma\left(  i\right)  $ and $\sigma\left(  j\right)
$. Hence, $\min\left(  \sigma\left(  \rho\left(  \left(  i,j\right)  \right)
\right)  \right)  =\sigma\left(  i\right)  $ and $\max\left(  \sigma\left(
\rho\left(  \left(  i,j\right)  \right)  \right)  \right)  =\sigma\left(
j\right)  $.
\par
Now, the definition of $a_{\sigma\left(  \rho\left(  \left(  i,j\right)
\right)  \right)  }$ shows that%
\begin{align*}
a_{\sigma\left(  \rho\left(  \left(  i,j\right)  \right)  \right)  }  &
=a_{\left(  \min\left(  \sigma\left(  \rho\left(  \left(  i,j\right)  \right)
\right)  \right)  ,\max\left(  \sigma\left(  \rho\left(  \left(  i,j\right)
\right)  \right)  \right)  \right)  }=a_{\left(  \sigma\left(  i\right)
,\sigma\left(  j\right)  \right)  }\\
&  \ \ \ \ \ \ \ \ \ \ \left(  \text{since }\min\left(  \sigma\left(
\rho\left(  \left(  i,j\right)  \right)  \right)  \right)  =\sigma\left(
i\right)  \text{ and }\max\left(  \sigma\left(  \rho\left(  \left(
i,j\right)  \right)  \right)  \right)  =\sigma\left(  j\right)  \right)  .
\end{align*}
This proves (\ref{sol.perm.sign.pseudoexplicit.short.b.a3}).}.
\end{itemize}

Now, recall that $G=\left\{  \left(  u,v\right)  \in\left[  n\right]
^{2}\ \mid\ u<v\right\}  $. Hence, the product sign $\prod_{1\leq i<j\leq n}$
in $\prod_{1\leq i<j\leq n}a_{\left(  \sigma\left(  i\right)  ,\sigma\left(
j\right)  \right)  }$ can be replaced by the product sign $\prod_{\left(
i,j\right)  \in G}$. Thus, we have%
\begin{align}
&  \underbrace{\prod_{1\leq i<j\leq n}}_{=\prod_{\left(  i,j\right)  \in G}%
}a_{\left(  \sigma\left(  i\right)  ,\sigma\left(  j\right)  \right)
}\nonumber\\
&  =\prod_{\left(  i,j\right)  \in G}a_{\left(  \sigma\left(  i\right)
,\sigma\left(  j\right)  \right)  }=\left(  \prod_{\substack{\left(
i,j\right)  \in G;\\\left(  i,j\right)  \in\operatorname*{Inv}\left(
\sigma\right)  }}\underbrace{a_{\left(  \sigma\left(  i\right)  ,\sigma\left(
j\right)  \right)  }}_{\substack{=-a_{\sigma\left(  \rho\left(  \left(
i,j\right)  \right)  \right)  }\\\text{(by
(\ref{sol.perm.sign.pseudoexplicit.short.b.a2}))}}}\right)  \cdot\left(
\prod_{\substack{\left(  i,j\right)  \in G;\\\left(  i,j\right)
\notin\operatorname*{Inv}\left(  \sigma\right)  }}\underbrace{a_{\left(
\sigma\left(  i\right)  ,\sigma\left(  j\right)  \right)  }}%
_{\substack{=a_{\sigma\left(  \rho\left(  \left(  i,j\right)  \right)
\right)  }\\\text{(by (\ref{sol.perm.sign.pseudoexplicit.short.b.a3}))}%
}}\right) \nonumber\\
&  =\underbrace{\left(  \prod_{\substack{\left(  i,j\right)  \in G;\\\left(
i,j\right)  \in\operatorname*{Inv}\left(  \sigma\right)  }}\left(
-a_{\sigma\left(  \rho\left(  \left(  i,j\right)  \right)  \right)  }\right)
\right)  }_{=\left(  -1\right)  ^{\left\vert \left\{  \left(  i,j\right)  \in
G\ \mid\ \left(  i,j\right)  \in\operatorname*{Inv}\left(  \sigma\right)
\right\}  \right\vert }\left(  \prod_{\substack{\left(  i,j\right)  \in
G;\\\left(  i,j\right)  \in\operatorname*{Inv}\left(  \sigma\right)
}}a_{\sigma\left(  \rho\left(  \left(  i,j\right)  \right)  \right)  }\right)
}\cdot\left(  \prod_{\substack{\left(  i,j\right)  \in G;\\\left(  i,j\right)
\notin\operatorname*{Inv}\left(  \sigma\right)  }}a_{\sigma\left(  \rho\left(
\left(  i,j\right)  \right)  \right)  }\right) \nonumber\\
&  =\underbrace{\left(  -1\right)  ^{\left\vert \left\{  \left(  i,j\right)
\in G\ \mid\ \left(  i,j\right)  \in\operatorname*{Inv}\left(  \sigma\right)
\right\}  \right\vert }}_{\substack{=\left(  -1\right)  ^{\left\vert
\operatorname*{Inv}\left(  \sigma\right)  \right\vert }\\\text{(since
}\operatorname*{Inv}\left(  \sigma\right)  \subseteq G\text{, and
thus}\\\left\{  \left(  i,j\right)  \in G\ \mid\ \left(  i,j\right)
\in\operatorname*{Inv}\left(  \sigma\right)  \right\}  =\operatorname*{Inv}%
\left(  \sigma\right)  \text{)}}}\underbrace{\left(  \prod_{\substack{\left(
i,j\right)  \in G;\\\left(  i,j\right)  \in\operatorname*{Inv}\left(
\sigma\right)  }}a_{\sigma\left(  \rho\left(  \left(  i,j\right)  \right)
\right)  }\right)  \cdot\left(  \prod_{\substack{\left(  i,j\right)  \in
G;\\\left(  i,j\right)  \notin\operatorname*{Inv}\left(  \sigma\right)
}}a_{\sigma\left(  \rho\left(  \left(  i,j\right)  \right)  \right)  }\right)
}_{=\prod_{\left(  i,j\right)  \in G}a_{\sigma\left(  \rho\left(  \left(
i,j\right)  \right)  \right)  }}\nonumber\\
&  =\underbrace{\left(  -1\right)  ^{\left\vert \operatorname*{Inv}\left(
\sigma\right)  \right\vert }}_{=\left(  -1\right)  ^{\sigma}}\underbrace{\prod
_{\left(  i,j\right)  \in G}a_{\sigma\left(  \rho\left(  \left(  i,j\right)
\right)  \right)  }}_{\substack{=\prod_{T\in\mathcal{P}_{2}\left(  \left[
n\right]  \right)  }a_{\sigma\left(  T\right)  }\\\text{(here, we have
substituted }T\\\text{for }\rho\left(  \left(  i,j\right)  \right)  \text{ in
the product,}\\\text{since the map }\rho:G\rightarrow\mathcal{P}_{2}\left(
\left[  n\right]  \right)  \\\text{is bijective)}}}=\left(  -1\right)
^{\sigma}\cdot\prod_{T\in\mathcal{P}_{2}\left(  \left[  n\right]  \right)
}a_{\sigma\left(  T\right)  }. \label{sol.perm.sign.pseudoexplicit.short.b.u}%
\end{align}
On the other hand,
\begin{align}
\underbrace{\prod_{1\leq i<j\leq n}}_{=\prod_{\left(  i,j\right)  \in G}%
}\ \ \underbrace{a_{\left(  i,j\right)  }}_{\substack{=a_{\rho\left(  \left(
i,j\right)  \right)  }\\\text{(by
(\ref{sol.perm.sign.pseudoexplicit.short.b.a1}))}}}  &  =\prod_{\left(
i,j\right)  \in G}a_{\rho\left(  \left(  i,j\right)  \right)  }=\prod
_{S\in\mathcal{P}_{2}\left(  \left[  n\right]  \right)  }a_{S}\nonumber\\
&  \ \ \ \ \ \ \ \ \ \ \left(
\begin{array}
[c]{c}%
\text{here, we have substituted }S\text{ for }\rho\left(  \left(  i,j\right)
\right)  \text{ in the product,}\\
\text{since the map }\rho:G\rightarrow\mathcal{P}_{2}\left(  \left[  n\right]
\right)  \text{ is bijective}%
\end{array}
\right) \nonumber\\
&  =\prod_{T\in\mathcal{P}_{2}\left(  \left[  n\right]  \right)
}\underbrace{a_{\sigma_{\ast}\left(  T\right)  }}_{\substack{=a_{\sigma\left(
T\right)  }\\\text{(since }\sigma_{\ast}\left(  T\right)  =\sigma\left(
T\right)  \\\text{(by the definition of }\sigma_{\ast}\text{))}}}\nonumber\\
&  \ \ \ \ \ \ \ \ \ \ \left(
\begin{array}
[c]{c}%
\text{here, we have substituted }\sigma_{\ast}\left(  T\right)  \text{ for
}S\text{ in the product,}\\
\text{since the map }\sigma_{\ast}:\mathcal{P}_{2}\left(  \left[  n\right]
\right)  \rightarrow\mathcal{P}_{2}\left(  \left[  n\right]  \right)  \text{
is bijective}%
\end{array}
\right) \nonumber\\
&  =\prod_{T\in\mathcal{P}_{2}\left(  \left[  n\right]  \right)  }%
a_{\sigma\left(  T\right)  }. \label{sol.perm.sign.pseudoexplicit.short.b.v}%
\end{align}
Hence, (\ref{sol.perm.sign.pseudoexplicit.short.b.u}) becomes%
\[
\prod_{1\leq i<j\leq n}a_{\left(  \sigma\left(  i\right)  ,\sigma\left(
j\right)  \right)  }=\left(  -1\right)  ^{\sigma}\cdot\underbrace{\prod
_{T\in\mathcal{P}_{2}\left(  \left[  n\right]  \right)  }a_{\sigma\left(
T\right)  }}_{\substack{=\prod_{1\leq i<j\leq n}a_{\left(  i,j\right)
}\\\text{(by (\ref{sol.perm.sign.pseudoexplicit.short.b.v}))}}}=\left(
-1\right)  ^{\sigma}\cdot\prod_{1\leq i<j\leq n}a_{\left(  i,j\right)  }.
\]
This solves Exercise \ref{exe.perm.sign.pseudoexplicit} \textbf{(b)}.

\textbf{(a)} We have $x_{j}-x_{i}=-\left(  x_{i}-x_{j}\right)  $ for every
$\left(  i,j\right)  \in\left\{  1,2,\ldots,n\right\}  ^{2}$. Hence, we can
apply Exercise \ref{exe.perm.sign.pseudoexplicit} \textbf{(b)} to $a_{\left(
i,j\right)  }=x_{i}-x_{j}$. As a result, we obtain $\prod_{1\leq i<j\leq
n}\left(  x_{\sigma\left(  i\right)  }-x_{\sigma\left(  j\right)  }\right)
=\left(  -1\right)  ^{\sigma}\cdot\prod_{1\leq i<j\leq n}\left(  x_{i}%
-x_{j}\right)  $. This solves Exercise \ref{exe.perm.sign.pseudoexplicit}
\textbf{(a)}.

\textbf{(c)} Applying Exercise \ref{exe.perm.sign.pseudoexplicit} \textbf{(a)}
to $x_{i}=i$, we obtain $\prod_{1\leq i<j\leq n}\left(  \sigma\left(
i\right)  -\sigma\left(  j\right)  \right)  =\left(  -1\right)  ^{\sigma}%
\cdot\prod_{1\leq i<j\leq n}\left(  i-j\right)  $. We can divide both sides of
this equality by $\prod_{1\leq i<j\leq n}\left(  i-j\right)  $ (because
$\prod_{1\leq i<j\leq n}\left(  i-j\right)  $ is a product of nonzero
integers, and thus nonzero). As a result, we obtain $\dfrac{\prod_{1\leq
i<j\leq n}\left(  \sigma\left(  i\right)  -\sigma\left(  j\right)  \right)
}{\prod_{1\leq i<j\leq n}\left(  i-j\right)  }=\left(  -1\right)  ^{\sigma}$.
Thus,%
\[
\left(  -1\right)  ^{\sigma}=\dfrac{\prod_{1\leq i<j\leq n}\left(
\sigma\left(  i\right)  -\sigma\left(  j\right)  \right)  }{\prod_{1\leq
i<j\leq n}\left(  i-j\right)  }=\prod_{1\leq i<j\leq n}\dfrac{\sigma\left(
i\right)  -\sigma\left(  j\right)  }{i-j}.
\]
Thus, (\ref{eq.sign.pseudoexplicit}) is proven. This solves Exercise
\ref{exe.perm.sign.pseudoexplicit} \textbf{(c)}.

\textbf{(d)} See below.
\end{proof}
\end{vershort}

\begin{verlong}
\begin{proof}
[Solution to Exercise \ref{exe.perm.sign.pseudoexplicit}.]\textbf{(b)} Let
$\left[  n\right]  $ be the set $\left\{  1,2,\ldots,n\right\}  $. Recall that
$S_{n}$ is the set of all permutations of the set $\left\{  1,2,\ldots
,n\right\}  $. In other words, $S_{n}$ is the set of all permutations of the
set $\left[  n\right]  $ (since $\left\{  1,2,\ldots,n\right\}  =\left[
n\right]  $).

Let $G$ be the subset%
\[
\left\{  \left(  i,j\right)  \in\left[  n\right]  ^{2}\ \mid\ i<j\right\}
\]
of $\left[  n\right]  ^{2}$. Thus,%
\[
G=\left\{  \left(  i,j\right)  \in\left[  n\right]  ^{2}\ \mid\ i<j\right\}
=\left\{  \left(  u,v\right)  \in\left[  n\right]  ^{2}\ \mid\ u<v\right\}
\]
(here, we renamed the index $\left(  i,j\right)  $ as $\left(  u,v\right)  $).
Clearly, the set $G$ is finite (since it is a subset of the finite set
$\left[  n\right]  ^{2}$).

Let $\operatorname*{Inv}\left(  \sigma\right)  $ be the set of inversions of
$\sigma$. Then, $\left(  -1\right)  ^{\sigma}=\left(  -1\right)  ^{\left\vert
\operatorname*{Inv}\left(  \sigma\right)  \right\vert }$%
\ \ \ \ \footnote{\textit{Proof.} We have%
\begin{align*}
\ell\left(  \sigma\right)   &  =\left(  \text{the number of inversions of
}\sigma\right)  \ \ \ \ \ \ \ \ \ \ \left(  \text{by the definition of }%
\ell\left(  \sigma\right)  \right) \\
&  =\left(  \text{the number of elements of }\operatorname*{Inv}\left(
\sigma\right)  \right) \\
&  \ \ \ \ \ \ \ \ \ \ \left(  \text{since }\operatorname*{Inv}\left(
\sigma\right)  \text{ is the set of all inversions of }\sigma\right) \\
&  =\left\vert \operatorname*{Inv}\left(  \sigma\right)  \right\vert .
\end{align*}
\par
But the definition of $\left(  -1\right)  ^{\sigma}$ yields $\left(
-1\right)  ^{\sigma}=\left(  -1\right)  ^{\ell\left(  \sigma\right)  }=\left(
-1\right)  ^{\left\vert \operatorname*{Inv}\left(  \sigma\right)  \right\vert
}$ (since $\ell\left(  \sigma\right)  =\left\vert \operatorname*{Inv}\left(
\sigma\right)  \right\vert $), qed.}. We notice furthermore that
$\operatorname*{Inv}\left(  \sigma\right)  \subseteq G$%
\ \ \ \ \footnote{\textit{Proof.} Let $c\in\operatorname*{Inv}\left(
\sigma\right)  $. Thus, $c$ is an inversion of $\sigma$ (since
$\operatorname*{Inv}\left(  \sigma\right)  $ is the set of inversions of
$\sigma$). In other words, $c$ is a pair $\left(  i,j\right)  $ of integers
satisfying $1\leq i<j\leq n$ and $\sigma\left(  i\right)  >\sigma\left(
j\right)  $. Consider this $\left(  i,j\right)  $. We have $i\in\left[
n\right]  $ (since $1\leq i\leq n$) and $j\in\left[  n\right]  $ (since $1\leq
j\leq n$). Thus, $\left(  i,j\right)  \in\left[  n\right]  ^{2}$ (since
$i\in\left[  n\right]  $ and $j\in\left[  n\right]  $). Hence, $\left(
i,j\right)  $ is an element of $\left[  n\right]  ^{2}$ and satisfies $i<j$.
In other words, $\left(  i,j\right)  $ is an element $\left(  u,v\right)  $ of
$\left[  n\right]  ^{2}$ and satisfies $u<v$. In other words, $\left(
i,j\right)  \in\left\{  \left(  u,v\right)  \in\left[  n\right]  ^{2}%
\ \mid\ u<v\right\}  =G$. Hence, $c=\left(  i,j\right)  \in G$.
\par
Now, let us forget that we fixed $c$. We thus have proven that every
$c\in\operatorname*{Inv}\left(  \sigma\right)  $ satisfies $c\in G$. In other
words, $\operatorname*{Inv}\left(  \sigma\right)  \subseteq G$, qed.}.

We shall use the Iverson bracket notation introduced in Definition
\ref{def.iverson}.

A useful property of the Iverson bracket is that it turns cardinalities of
sets into sums: Namely, if $S$ is a finite set, and if $T$ is a subset of $S$,
then%
\begin{equation}
\left\vert T\right\vert =\sum_{s\in S}\left[  s\in T\right]
\label{sol.perm.sign.pseudoexplicit.b.iverson.sum}%
\end{equation}
\footnote{In fact, this is precisely the claim of Lemma \ref{lem.iverson.card}%
.}. Consequently, if $S$ is a finite set, and if $T$ is a subset of $S$, then%
\begin{equation}
\left(  -1\right)  ^{\left\vert T\right\vert }=\prod_{s\in S}\left(
-1\right)  ^{\left[  s\in T\right]  }
\label{sol.perm.sign.pseudoexplicit.b.iverson.prod}%
\end{equation}
\footnote{\textit{Proof of (\ref{sol.perm.sign.pseudoexplicit.b.iverson.prod}%
):} Let $S$ be a finite set. Let $T$ be a subset of $S$. Then,
(\ref{sol.perm.sign.pseudoexplicit.b.iverson.sum}) yields $\left\vert
T\right\vert =\sum_{s\in S}\left[  s\in T\right]  $. Hence,
\[
\left(  -1\right)  ^{\left\vert T\right\vert }=\left(  -1\right)  ^{\sum_{s\in
S}\left[  s\in T\right]  }=\prod_{s\in S}\left(  -1\right)  ^{\left[  s\in
T\right]  }%
\]
(because $a^{\sum_{s\in S}b_{s}}=\prod_{s\in S}a^{b_{s}}$ whenever $a$ is an
integer and $b_{s}$ is an integer for every $s\in S$). This proves
(\ref{sol.perm.sign.pseudoexplicit.b.iverson.prod}).}. We can apply this to
$S=G$ and $T=\operatorname*{Inv}\left(  \sigma\right)  $ (since
$\operatorname*{Inv}\left(  \sigma\right)  \subseteq G$). As a result, we
obtain%
\[
\left(  -1\right)  ^{\left\vert \operatorname*{Inv}\left(  \sigma\right)
\right\vert }=\prod_{s\in G}\left(  -1\right)  ^{\left[  s\in
\operatorname*{Inv}\left(  \sigma\right)  \right]  }=\prod_{\left(
i,j\right)  \in G}\left(  -1\right)  ^{\left[  \left(  i,j\right)
\in\operatorname*{Inv}\left(  \sigma\right)  \right]  }%
\]
(here, we renamed the index $s$ as $\left(  i,j\right)  $, because all
elements of $G$ are pairs). Thus,%
\begin{equation}
\left(  -1\right)  ^{\sigma}=\left(  -1\right)  ^{\left\vert
\operatorname*{Inv}\left(  \sigma\right)  \right\vert }=\prod_{\left(
i,j\right)  \in G}\left(  -1\right)  ^{\left[  \left(  i,j\right)
\in\operatorname*{Inv}\left(  \sigma\right)  \right]  }.
\label{sol.perm.sign.pseudoexplicit.b.2}%
\end{equation}

On the other hand, for every $\tau\in S_{n}$, we define a map $\tau^{\left[
2\right]  }:G\rightarrow G$ by setting%
\begin{equation}
\left(  \tau^{\left[  2\right]  }\left(  i,j\right)  =\left(  \min\left\{
\tau\left(  i\right)  ,\tau\left(  j\right)  \right\}  ,\max\left\{
\tau\left(  i\right)  ,\tau\left(  j\right)  \right\}  \right)
\ \ \ \ \ \ \ \ \ \ \text{for every }\left(  i,j\right)  \in G\right)  .
\label{sol.perm.sign.pseudoexplicit.b.tau}%
\end{equation}
This map $\tau^{\left[  2\right]  }$ is well-defined\footnote{\textit{Proof.}
In order to prove this, we need to show that every element of $G$ can be
written in the form $\left(  i,j\right)  $, and that every $\left(
i,j\right)  \in G$ satisfies $\left(  \min\left\{  \tau\left(  i\right)
,\tau\left(  j\right)  \right\}  ,\max\left\{  \tau\left(  i\right)
,\tau\left(  j\right)  \right\}  \right)  \in G$. The first of these two
claims (i.e., that every element of $G$ can be written in the form $\left(
i,j\right)  $) is obvious. It thus remains to prove the second of these two
claims. In other words, it remains to prove that every $\left(  i,j\right)
\in G$ satisfies $\left(  \min\left\{  \tau\left(  i\right)  ,\tau\left(
j\right)  \right\}  ,\max\left\{  \tau\left(  i\right)  ,\tau\left(  j\right)
\right\}  \right)  \in G$.
\par
So let $\left(  i,j\right)  \in G$. We must prove that $\left(  \min\left\{
\tau\left(  i\right)  ,\tau\left(  j\right)  \right\}  ,\max\left\{
\tau\left(  i\right)  ,\tau\left(  j\right)  \right\}  \right)  \in G$.
\par
We have $\left(  i,j\right)  \in G=\left\{  \left(  u,v\right)  \in\left[
n\right]  ^{2}\ \mid\ u<v\right\}  $. In other words, $\left(  i,j\right)  $
is an element $\left(  u,v\right)  $ of $\left[  n\right]  ^{2}$ satisfying
$u<v$. In other words, $\left(  i,j\right)  $ is an element of $\left[
n\right]  ^{2}$ satisfying $i<j$.
\par
Now, $\tau\in S_{n}$. Hence, $\tau$ is a permutation of the set $\left[
n\right]  $ (since $S_{n}$ is the set of all permutations of the set $\left[
n\right]  $). In other words, $\tau$ is a bijective map $\left[  n\right]
\rightarrow\left[  n\right]  $. Hence, $\tau$ is both surjective and
injective. We have $i\neq j$ (since $i<j$) and thus $\tau\left(  i\right)
\neq\tau\left(  j\right)  $ (since $\tau$ is injective). Clearly, both
$\tau\left(  i\right)  $ and $\tau\left(  j\right)  $ are elements of $\left[
n\right]  $. Thus, $\left\{  \tau\left(  i\right)  ,\tau\left(  j\right)
\right\}  \subseteq\left[  n\right]  $, so that $\min\left\{  \tau\left(
i\right)  ,\tau\left(  j\right)  \right\}  \in\left\{  \tau\left(  i\right)
,\tau\left(  j\right)  \right\}  \subseteq\left[  n\right]  $ and
$\max\left\{  \tau\left(  i\right)  ,\tau\left(  j\right)  \right\}
\in\left\{  \tau\left(  i\right)  ,\tau\left(  j\right)  \right\}
\subseteq\left[  n\right]  $. Hence,
\[
\left(  \underbrace{\min\left\{  \tau\left(  i\right)  ,\tau\left(  j\right)
\right\}  }_{\in\left[  n\right]  },\underbrace{\max\left\{  \tau\left(
i\right)  ,\tau\left(  j\right)  \right\}  }_{\in\left[  n\right]  }\right)
\in\left[  n\right]  \times\left[  n\right]  =\left[  n\right]  ^{2}.
\]
\par
Now, let us show that $\min\left\{  \tau\left(  i\right)  ,\tau\left(
j\right)  \right\}  <\max\left\{  \tau\left(  i\right)  ,\tau\left(  j\right)
\right\}  $. Indeed, we assume the contrary (for the sake of contradiction).
Thus, $\min\left\{  \tau\left(  i\right)  ,\tau\left(  j\right)  \right\}
\geq\max\left\{  \tau\left(  i\right)  ,\tau\left(  j\right)  \right\}  $.
\par
Now, $\tau\left(  i\right)  $ is an element of the set $\left\{  \tau\left(
i\right)  ,\tau\left(  j\right)  \right\}  $, and thus greater or equal to the
minimum of this set. In other words, we have $\tau\left(  i\right)  \geq
\min\left\{  \tau\left(  i\right)  ,\tau\left(  j\right)  \right\}  $.
\par
Also, $\tau\left(  i\right)  $ is an element of the set $\left\{  \tau\left(
i\right)  ,\tau\left(  j\right)  \right\}  $, and thus less or equal to the
maximum of this set. In other words, we have $\tau\left(  i\right)  \leq
\max\left\{  \tau\left(  i\right)  ,\tau\left(  j\right)  \right\}  $.
\par
Combining $\tau\left(  i\right)  \geq\min\left\{  \tau\left(  i\right)
,\tau\left(  j\right)  \right\}  \geq\max\left\{  \tau\left(  i\right)
,\tau\left(  j\right)  \right\}  $ with $\tau\left(  i\right)  \leq
\max\left\{  \tau\left(  i\right)  ,\tau\left(  j\right)  \right\}  $, we
obtain $\tau\left(  i\right)  =\max\left\{  \tau\left(  i\right)  ,\tau\left(
j\right)  \right\}  $.
\par
Also, $\tau\left(  j\right)  $ is an element of the set $\left\{  \tau\left(
i\right)  ,\tau\left(  j\right)  \right\}  $, and thus greater or equal to the
minimum of this set. In other words, we have $\tau\left(  j\right)  \geq
\min\left\{  \tau\left(  i\right)  ,\tau\left(  j\right)  \right\}  $.
\par
Also, $\tau\left(  j\right)  $ is an element of the set $\left\{  \tau\left(
i\right)  ,\tau\left(  j\right)  \right\}  $, and thus less or equal to the
maximum of this set. In other words, we have $\tau\left(  j\right)  \leq
\max\left\{  \tau\left(  i\right)  ,\tau\left(  j\right)  \right\}  $.
\par
Combining $\tau\left(  j\right)  \geq\min\left\{  \tau\left(  i\right)
,\tau\left(  j\right)  \right\}  \geq\max\left\{  \tau\left(  i\right)
,\tau\left(  j\right)  \right\}  $ with $\tau\left(  j\right)  \leq
\max\left\{  \tau\left(  i\right)  ,\tau\left(  j\right)  \right\}  $, we
obtain $\tau\left(  j\right)  =\max\left\{  \tau\left(  i\right)  ,\tau\left(
j\right)  \right\}  $. Compared with $\tau\left(  i\right)  =\max\left\{
\tau\left(  i\right)  ,\tau\left(  j\right)  \right\}  $, this yields
$\tau\left(  i\right)  =\tau\left(  j\right)  $, which contradicts
$\tau\left(  i\right)  \neq\tau\left(  j\right)  $. This contradiction proves
that our assumption was wrong.
\par
Hence, $\min\left\{  \tau\left(  i\right)  ,\tau\left(  j\right)  \right\}
<\max\left\{  \tau\left(  i\right)  ,\tau\left(  j\right)  \right\}  $ is
proven. Now, we know that $\left(  \min\left\{  \tau\left(  i\right)
,\tau\left(  j\right)  \right\}  ,\max\left\{  \tau\left(  i\right)
,\tau\left(  j\right)  \right\}  \right)  $ is an element of $\left[
n\right]  ^{2}$ satisfying $\min\left\{  \tau\left(  i\right)  ,\tau\left(
j\right)  \right\}  <\max\left\{  \tau\left(  i\right)  ,\tau\left(  j\right)
\right\}  $. In other words, $\left(  \min\left\{  \tau\left(  i\right)
,\tau\left(  j\right)  \right\}  ,\max\left\{  \tau\left(  i\right)
,\tau\left(  j\right)  \right\}  \right)  $ is an element $\left(  u,v\right)
\in\left[  n\right]  ^{2}$ satisfying $u<v$. In other words,%
\[
\left(  \min\left\{  \tau\left(  i\right)  ,\tau\left(  j\right)  \right\}
,\max\left\{  \tau\left(  i\right)  ,\tau\left(  j\right)  \right\}  \right)
\in\left\{  \left(  u,v\right)  \in\left[  n\right]  ^{2}\ \mid\ u<v\right\}
=G.
\]
Thus, we have proven that $\left(  \min\left\{  \tau\left(  i\right)
,\tau\left(  j\right)  \right\}  ,\max\left\{  \tau\left(  i\right)
,\tau\left(  j\right)  \right\}  \right)  \in G$. Qed.}. Moreover, we have%
\begin{equation}
\left(  \tau^{-1}\right)  ^{\left[  2\right]  }\circ\tau^{\left[  2\right]
}=\operatorname*{id}\ \ \ \ \ \ \ \ \ \ \text{for every }\tau\in S_{n}
\label{sol.perm.sign.pseudoexplicit.b.tautau}%
\end{equation}
\footnote{\textit{Proof of (\ref{sol.perm.sign.pseudoexplicit.b.tautau}):} Let
$\tau\in S_{n}$. Let $c\in G$. We are going to prove that $\left(  \left(
\tau^{-1}\right)  ^{\left[  2\right]  }\circ\tau^{\left[  2\right]  }\right)
\left(  c\right)  =c$.
\par
We have $c\in G=\left\{  \left(  i,j\right)  \in\left[  n\right]  ^{2}%
\ \mid\ i<j\right\}  $. In other words, $c$ can be written in the form
$c=\left(  i,j\right)  $ for some $\left(  i,j\right)  \in\left[  n\right]
^{2}$ satisfying $i<j$. Consider this $\left(  i,j\right)  $.
\par
The definition of $\tau^{\left[  2\right]  }$ yields $\tau^{\left[  2\right]
}\left(  i,j\right)  =\left(  \min\left\{  \tau\left(  i\right)  ,\tau\left(
j\right)  \right\}  ,\max\left\{  \tau\left(  i\right)  ,\tau\left(  j\right)
\right\}  \right)  $. Thus, $\tau^{\left[  2\right]  }\left(  \underbrace{c}%
_{=\left(  i,j\right)  }\right)  =\tau^{\left[  2\right]  }\left(  i,j\right)
=\left(  \min\left\{  \tau\left(  i\right)  ,\tau\left(  j\right)  \right\}
,\max\left\{  \tau\left(  i\right)  ,\tau\left(  j\right)  \right\}  \right)
$.
\par
Now, $\tau\in S_{n}$. Hence, $\tau$ is a permutation of the set $\left[
n\right]  $ (since $S_{n}$ is the set of all permutations of the set $\left[
n\right]  $). In other words, $\tau$ is a bijective map $\left[  n\right]
\rightarrow\left[  n\right]  $. Hence, $\tau$ is both surjective and
injective. We have $i\neq j$ (since $i<j$) and thus $\tau\left(  i\right)
\neq\tau\left(  j\right)  $ (since $\tau$ is injective). Hence, we have either
$\tau\left(  i\right)  <\tau\left(  j\right)  $ or $\tau\left(  i\right)
>\tau\left(  j\right)  $. In other words, we are in one of the following two
cases:
\par
\textit{Case 1:} We have $\tau\left(  i\right)  <\tau\left(  j\right)  $.
\par
\textit{Case 2:} We have $\tau\left(  i\right)  >\tau\left(  j\right)  $.
\par
Let us first consider Case 1. In this case, we have $\tau\left(  i\right)
<\tau\left(  j\right)  $. Thus, $\min\left\{  \tau\left(  i\right)
,\tau\left(  j\right)  \right\}  =\tau\left(  i\right)  $ and $\max\left\{
\tau\left(  i\right)  ,\tau\left(  j\right)  \right\}  =\tau\left(  j\right)
$. Hence,%
\[
\tau^{\left[  2\right]  }\left(  c\right)  =\left(  \underbrace{\min\left\{
\tau\left(  i\right)  ,\tau\left(  j\right)  \right\}  }_{=\tau\left(
i\right)  },\underbrace{\max\left\{  \tau\left(  i\right)  ,\tau\left(
j\right)  \right\}  }_{=\tau\left(  j\right)  }\right)  =\left(  \tau\left(
i\right)  ,\tau\left(  j\right)  \right)  .
\]
Now,%
\begin{align*}
\left(  \left(  \tau^{-1}\right)  ^{\left[  2\right]  }\circ\tau^{\left[
2\right]  }\right)  \left(  c\right)   &  =\left(  \tau^{-1}\right)  ^{\left[
2\right]  }\left(  \underbrace{\tau^{\left[  2\right]  }\left(  c\right)
}_{=\left(  \tau\left(  i\right)  ,\tau\left(  j\right)  \right)  }\right)
=\left(  \tau^{-1}\right)  ^{\left[  2\right]  }\left(  \tau\left(  i\right)
,\tau\left(  j\right)  \right) \\
&  =\left(  \min\left\{  \underbrace{\tau^{-1}\left(  \tau\left(  i\right)
\right)  }_{=i},\underbrace{\tau^{-1}\left(  \tau\left(  j\right)  \right)
}_{=j}\right\}  ,\max\left\{  \underbrace{\tau^{-1}\left(  \tau\left(
i\right)  \right)  }_{=i},\underbrace{\tau^{-1}\left(  \tau\left(  j\right)
\right)  }_{=j}\right\}  \right) \\
&  \ \ \ \ \ \ \ \ \ \ \left(  \text{by the definition of }\left(  \tau
^{-1}\right)  ^{\left[  2\right]  }\right) \\
&  =\left(  \underbrace{\min\left\{  i,j\right\}  }%
_{\substack{=i\\\text{(since }i<j\text{)}}},\underbrace{\max\left\{
i,j\right\}  }_{\substack{=j\\\text{(since }i<j\text{)}}}\right)  =\left(
i,j\right)  =c.
\end{align*}
Hence, $\left(  \left(  \tau^{-1}\right)  ^{\left[  2\right]  }\circ
\tau^{\left[  2\right]  }\right)  \left(  c\right)  =c$ is proven in Case 1.
\par
Let us now consider Case 2. In this case, we have $\tau\left(  i\right)
>\tau\left(  j\right)  $. Thus, $\min\left\{  \tau\left(  i\right)
,\tau\left(  j\right)  \right\}  =\tau\left(  j\right)  $ and $\max\left\{
\tau\left(  i\right)  ,\tau\left(  j\right)  \right\}  =\tau\left(  i\right)
$. Hence,%
\[
\tau^{\left[  2\right]  }\left(  c\right)  =\left(  \underbrace{\min\left\{
\tau\left(  i\right)  ,\tau\left(  j\right)  \right\}  }_{=\tau\left(
j\right)  },\underbrace{\max\left\{  \tau\left(  i\right)  ,\tau\left(
j\right)  \right\}  }_{=\tau\left(  i\right)  }\right)  =\left(  \tau\left(
j\right)  ,\tau\left(  i\right)  \right)  .
\]
Now,%
\begin{align*}
\left(  \left(  \tau^{-1}\right)  ^{\left[  2\right]  }\circ\tau^{\left[
2\right]  }\right)  \left(  c\right)   &  =\left(  \tau^{-1}\right)  ^{\left[
2\right]  }\left(  \underbrace{\tau^{\left[  2\right]  }\left(  c\right)
}_{=\left(  \tau\left(  j\right)  ,\tau\left(  i\right)  \right)  }\right)
=\left(  \tau^{-1}\right)  ^{\left[  2\right]  }\left(  \tau\left(  j\right)
,\tau\left(  i\right)  \right) \\
&  =\left(  \min\left\{  \underbrace{\tau^{-1}\left(  \tau\left(  j\right)
\right)  }_{=j},\underbrace{\tau^{-1}\left(  \tau\left(  i\right)  \right)
}_{=i}\right\}  ,\max\left\{  \underbrace{\tau^{-1}\left(  \tau\left(
j\right)  \right)  }_{=j},\underbrace{\tau^{-1}\left(  \tau\left(  i\right)
\right)  }_{=i}\right\}  \right) \\
&  \ \ \ \ \ \ \ \ \ \ \left(  \text{by the definition of }\left(  \tau
^{-1}\right)  ^{\left[  2\right]  }\right) \\
&  =\left(  \underbrace{\min\left\{  j,i\right\}  }%
_{\substack{=i\\\text{(since }i<j\text{)}}},\underbrace{\max\left\{
j,i\right\}  }_{\substack{=j\\\text{(since }i<j\text{)}}}\right)  =\left(
i,j\right)  =c.
\end{align*}
Hence, $\left(  \left(  \tau^{-1}\right)  ^{\left[  2\right]  }\circ
\tau^{\left[  2\right]  }\right)  \left(  c\right)  =c$ is proven in Case 2.
\par
We have now proven $\left(  \left(  \tau^{-1}\right)  ^{\left[  2\right]
}\circ\tau^{\left[  2\right]  }\right)  \left(  c\right)  =c$ in each of the
two Cases 1 and 2. Thus, $\left(  \left(  \tau^{-1}\right)  ^{\left[
2\right]  }\circ\tau^{\left[  2\right]  }\right)  \left(  c\right)  =c$ always
holds.
\par
So we have $\left(  \left(  \tau^{-1}\right)  ^{\left[  2\right]  }\circ
\tau^{\left[  2\right]  }\right)  \left(  c\right)  =c=\operatorname*{id}%
\left(  c\right)  $.
\par
Let us now forget that we fixed $c$. We thus have proven that $\left(  \left(
\tau^{-1}\right)  ^{\left[  2\right]  }\circ\tau^{\left[  2\right]  }\right)
\left(  c\right)  =\operatorname*{id}\left(  c\right)  $ for every $c\in G$.
In other words, $\left(  \tau^{-1}\right)  ^{\left[  2\right]  }\circ
\tau^{\left[  2\right]  }=\operatorname*{id}$. This proves
(\ref{sol.perm.sign.pseudoexplicit.b.tautau}).}. Thus, for every $\tau\in
S_{n}$, the map $\tau^{\left[  2\right]  }$ is a
bijection\footnote{\textit{Proof.} Let $\tau\in S_{n}$. Applying
(\ref{sol.perm.sign.pseudoexplicit.b.tautau}) to $\tau^{-1}$ instead of $\tau
$, we obtain $\left(  \left(  \tau^{-1}\right)  ^{-1}\right)  ^{\left[
2\right]  }\circ\left(  \tau^{-1}\right)  ^{\left[  2\right]  }%
=\operatorname*{id}$. Since $\left(  \tau^{-1}\right)  ^{-1}=\tau$, this
rewrites as $\tau^{\left[  2\right]  }\circ\left(  \tau^{-1}\right)  ^{\left[
2\right]  }=\operatorname*{id}$. Combined with $\left(  \tau^{-1}\right)
^{\left[  2\right]  }\circ\tau^{\left[  2\right]  }=\operatorname*{id}$ (which
follows from (\ref{sol.perm.sign.pseudoexplicit.b.tautau})), this yields that
the maps $\tau^{\left[  2\right]  }$ and $\left(  \tau^{-1}\right)  ^{\left[
2\right]  }$ are mutually inverse. Hence, the map $\tau^{\left[  2\right]  }$
is invertible, i.e., is a bijection. Qed.}. Applying this to $\tau=\sigma$, we
see that the map $\sigma^{\left[  2\right]  }$ is a bijection.

Now, every $\left(  i,j\right)  \in G$ satisfies%
\begin{equation}
\left(  -1\right)  ^{\left[  \left(  i,j\right)  \in\operatorname*{Inv}\left(
\sigma\right)  \right]  }a_{\sigma^{\left[  2\right]  }\left(  i,j\right)
}=a_{\left(  \sigma\left(  i\right)  ,\sigma\left(  j\right)  \right)  }
\label{sol.perm.sign.pseudoexplicit.b.a-inv}%
\end{equation}
\footnote{\textit{Proof of (\ref{sol.perm.sign.pseudoexplicit.b.a-inv}):} Let
$\left(  i,j\right)  \in G$. Thus, $\left(  i,j\right)  \in G=\left\{  \left(
u,v\right)  \in\left[  n\right]  ^{2}\ \mid\ u<v\right\}  $. In other words,
$\left(  i,j\right)  $ is an element $\left(  u,v\right)  $ of $\left[
n\right]  ^{2}$ satisfying $u<v$. In other words, $\left(  i,j\right)  $ is an
element of $\left[  n\right]  ^{2}$ satisfying $i<j$. From $\left(
i,j\right)  \in\left[  n\right]  ^{2}$, we obtain $i\in\left[  n\right]  $ and
$j\in\left[  n\right]  $. Thus, $1\leq i$ (since $i\in\left[  n\right]  $) and
$j\leq n$ (since $j\in\left[  n\right]  $), so that $1\leq i<j\leq n$.
\par
We must be in one of the following two cases:
\par
\textit{Case 1:} We have $\sigma\left(  i\right)  \leq\sigma\left(  j\right)
$.
\par
\textit{Case 2:} We have $\sigma\left(  i\right)  >\sigma\left(  j\right)  $.
\par
Let us first consider Case 1. In this case, we have $\sigma\left(  i\right)
\leq\sigma\left(  j\right)  $. Thus, $\min\left\{  \sigma\left(  i\right)
,\sigma\left(  j\right)  \right\}  =\sigma\left(  i\right)  $ and
$\max\left\{  \sigma\left(  i\right)  ,\sigma\left(  j\right)  \right\}
=\sigma\left(  j\right)  $. Now, the definition of $\sigma^{\left[  2\right]
}$ yields%
\[
\sigma^{\left[  2\right]  }\left(  i,j\right)  =\left(  \underbrace{\min
\left\{  \sigma\left(  i\right)  ,\sigma\left(  j\right)  \right\}  }%
_{=\sigma\left(  i\right)  },\underbrace{\max\left\{  \sigma\left(  i\right)
,\sigma\left(  j\right)  \right\}  }_{=\sigma\left(  j\right)  }\right)
=\left(  \sigma\left(  i\right)  ,\sigma\left(  j\right)  \right)  .
\]
Hence, $a_{\sigma^{\left[  2\right]  }\left(  i,j\right)  }=a_{\left(
\sigma\left(  i\right)  ,\sigma\left(  j\right)  \right)  }$.
\par
On the other hand, let us show that $\left(  i,j\right)  \notin%
\operatorname*{Inv}\left(  \sigma\right)  $. Indeed, we assume the contrary
(for the sake of contradiction). Thus, $\left(  i,j\right)  \in
\operatorname*{Inv}\left(  \sigma\right)  $. In other words, $\left(
i,j\right)  $ is an inversion of $\sigma$ (since $\operatorname*{Inv}\left(
\sigma\right)  $ is the set of inversions of $\sigma$). In other words,
$\left(  i,j\right)  $ is a pair of integers satisfying $1\leq i<j\leq n$ and
$\sigma\left(  i\right)  >\sigma\left(  j\right)  $ (because of how an
\textquotedblleft inversion\textquotedblright\ is defined). But $\sigma\left(
i\right)  >\sigma\left(  j\right)  $ contradicts $\sigma\left(  i\right)
\leq\sigma\left(  j\right)  $. Thus, we have obtained a contradiction.
Therefore, our assumption must have been false. This proves that $\left(
i,j\right)  \notin\operatorname*{Inv}\left(  \sigma\right)  $.
\par
Hence, $\left(  i,j\right)  \in\operatorname*{Inv}\left(  \sigma\right)  $ is
false. Thus, $\left[  \left(  i,j\right)  \in\operatorname*{Inv}\left(
\sigma\right)  \right]  =0$, so that $\left(  -1\right)  ^{\left[  \left(
i,j\right)  \in\operatorname*{Inv}\left(  \sigma\right)  \right]  }=\left(
-1\right)  ^{0}=1$ and therefore $\underbrace{\left(  -1\right)  ^{\left[
\left(  i,j\right)  \in\operatorname*{Inv}\left(  \sigma\right)  \right]  }%
}_{=1}\underbrace{a_{\sigma^{\left[  2\right]  }\left(  i,j\right)  }%
}_{=a_{\left(  \sigma\left(  i\right)  ,\sigma\left(  j\right)  \right)  }%
}=1a_{\left(  \sigma\left(  i\right)  ,\sigma\left(  j\right)  \right)
}=a_{\left(  \sigma\left(  i\right)  ,\sigma\left(  j\right)  \right)  }$.
Thus, (\ref{sol.perm.sign.pseudoexplicit.b.a-inv}) is proven in Case 1.
\par
Let us now consider Case 2. In this case, we have $\sigma\left(  i\right)
>\sigma\left(  j\right)  $. Thus, $\min\left\{  \sigma\left(  i\right)
,\sigma\left(  j\right)  \right\}  =\sigma\left(  j\right)  $ and
$\max\left\{  \sigma\left(  i\right)  ,\sigma\left(  j\right)  \right\}
=\sigma\left(  i\right)  $. Now, the definition of $\sigma^{\left[  2\right]
}$ yields%
\[
\sigma^{\left[  2\right]  }\left(  i,j\right)  =\left(  \underbrace{\min
\left\{  \sigma\left(  i\right)  ,\sigma\left(  j\right)  \right\}  }%
_{=\sigma\left(  j\right)  },\underbrace{\max\left\{  \sigma\left(  i\right)
,\sigma\left(  j\right)  \right\}  }_{=\sigma\left(  i\right)  }\right)
=\left(  \sigma\left(  j\right)  ,\sigma\left(  i\right)  \right)  .
\]
Hence, $a_{\sigma^{\left[  2\right]  }\left(  i,j\right)  }=a_{\left(
\sigma\left(  j\right)  ,\sigma\left(  i\right)  \right)  }=-a_{\left(
\sigma\left(  i\right)  ,\sigma\left(  j\right)  \right)  }$ (by
(\ref{eq.exe.perm.sign.pseudoexplicit.b.skew}), applied to $\sigma\left(
i\right)  $ and $\sigma\left(  j\right)  $ instead of $i$ and $j$).
\par
On the other hand, $\left(  i,j\right)  $ is a pair of integers satisfying
$1\leq i<j\leq n$ and $\sigma\left(  i\right)  >\sigma\left(  j\right)  $. In
other words, $\left(  i,j\right)  $ is an inversion of $\sigma$ (because of
how an \textquotedblleft inversion\textquotedblright\ is defined). In other
words, $\left(  i,j\right)  \in\operatorname*{Inv}\left(  \sigma\right)  $
(since $\operatorname*{Inv}\left(  \sigma\right)  $ is the set of inversions
of $\sigma$). Hence, $\left[  \left(  i,j\right)  \in\operatorname*{Inv}%
\left(  \sigma\right)  \right]  =1$, so that $\left(  -1\right)  ^{\left[
\left(  i,j\right)  \in\operatorname*{Inv}\left(  \sigma\right)  \right]
}=\left(  -1\right)  ^{1}=-1$ and therefore $\underbrace{\left(  -1\right)
^{\left[  \left(  i,j\right)  \in\operatorname*{Inv}\left(  \sigma\right)
\right]  }}_{=-1}\underbrace{a_{\sigma^{\left[  2\right]  }\left(  i,j\right)
}}_{=-a_{\left(  \sigma\left(  i\right)  ,\sigma\left(  j\right)  \right)  }%
}=\left(  -1\right)  \left(  -a_{\left(  \sigma\left(  i\right)
,\sigma\left(  j\right)  \right)  }\right)  =a_{\left(  \sigma\left(
i\right)  ,\sigma\left(  j\right)  \right)  }$. Thus,
(\ref{sol.perm.sign.pseudoexplicit.b.a-inv}) is proven in Case 2.
\par
We have now proven (\ref{sol.perm.sign.pseudoexplicit.b.a-inv}) in each of the
two Cases 1 and 2. Thus, (\ref{sol.perm.sign.pseudoexplicit.b.a-inv}) always
holds, qed.}.

Now,%
\begin{align*}
&  \underbrace{\prod_{1\leq i<j\leq n}}_{\substack{=\prod_{\substack{\left(
i,j\right)  \in\left[  n\right]  ^{2};\\i<j}}=\prod_{\left(  i,j\right)  \in
G}\\\text{(since }G=\left\{  \left(  i,j\right)  \in\left[  n\right]
^{2}\ \mid\ i<j\right\}  \text{)}}}\underbrace{a_{\left(  \sigma\left(
i\right)  ,\sigma\left(  j\right)  \right)  }}_{\substack{=\left(  -1\right)
^{\left[  \left(  i,j\right)  \in\operatorname*{Inv}\left(  \sigma\right)
\right]  }a_{\sigma^{\left[  2\right]  }\left(  i,j\right)  }\\\text{(by
(\ref{sol.perm.sign.pseudoexplicit.b.a-inv}))}}}\\
&  =\prod_{\left(  i,j\right)  \in G}\left(  \left(  -1\right)  ^{\left[
\left(  i,j\right)  \in\operatorname*{Inv}\left(  \sigma\right)  \right]
}a_{\sigma^{\left[  2\right]  }\left(  i,j\right)  }\right)
=\underbrace{\left(  \prod_{\left(  i,j\right)  \in G}\left(  -1\right)
^{\left[  \left(  i,j\right)  \in\operatorname*{Inv}\left(  \sigma\right)
\right]  }\right)  }_{\substack{=\left(  -1\right)  ^{\sigma}\\\text{(by
(\ref{sol.perm.sign.pseudoexplicit.b.2}))}}}\left(  \prod_{\left(  i,j\right)
\in G}a_{\sigma^{\left[  2\right]  }\left(  i,j\right)  }\right) \\
&  =\left(  -1\right)  ^{\sigma}\prod_{\left(  i,j\right)  \in G}%
a_{\sigma^{\left[  2\right]  }\left(  i,j\right)  }=\left(  -1\right)
^{\sigma}\prod_{c\in G}a_{\sigma^{\left[  2\right]  }\left(  c\right)  }\\
&  \ \ \ \ \ \ \ \ \ \ \left(  \text{here, we have renamed the index }\left(
i,j\right)  \text{ as }c\text{ in the product}\right) \\
&  =\left(  -1\right)  ^{\sigma}\prod_{c\in G}a_{c}\ \ \ \ \ \ \ \ \ \ \left(
\begin{array}
[c]{c}%
\text{here, we have substituted }c\text{ for }\sigma^{\left[  2\right]
}\left(  c\right)  \text{ in the product,}\\
\text{since the map }\sigma^{\left[  2\right]  }:G\rightarrow G\text{ is a
bijection}%
\end{array}
\right) \\
&  =\left(  -1\right)  ^{\sigma}\underbrace{\prod_{\left(  i,j\right)  \in G}%
}_{\substack{=\prod_{\substack{\left(  i,j\right)  \in\left[  n\right]
^{2};\\i<j}}\\\text{(since }G=\left\{  \left(  i,j\right)  \in\left[
n\right]  ^{2}\ \mid\ i<j\right\}  \text{)}}}a_{\left(  i,j\right)  }\\
&  \ \ \ \ \ \ \ \ \ \ \left(
\begin{array}
[c]{c}%
\text{here, we have renamed the index }c\text{ as }\left(  i,j\right)  \text{
in the product,}\\
\text{since every element of }G\text{ has the form }\left(  i,j\right)
\end{array}
\right) \\
&  =\left(  -1\right)  ^{\sigma}\underbrace{\prod_{\substack{\left(
i,j\right)  \in\left[  n\right]  ^{2};\\i<j}}}_{=\prod_{1\leq i<j\leq n}%
}a_{\left(  i,j\right)  }=\left(  -1\right)  ^{\sigma}\prod_{1\leq i<j\leq
n}a_{\left(  i,j\right)  }.
\end{align*}
This solves Exercise \ref{exe.perm.sign.pseudoexplicit} \textbf{(b)}.

\textbf{(a)} We have $x_{j}-x_{i}=-\left(  x_{i}-x_{j}\right)  $ for every
$\left(  i,j\right)  \in\left\{  1,2,\ldots,n\right\}  ^{2}$. Hence, we can
apply Exercise \ref{exe.perm.sign.pseudoexplicit} \textbf{(b)} to $a_{\left(
i,j\right)  }=x_{i}-x_{j}$. As a result, we obtain $\prod_{1\leq i<j\leq
n}\left(  x_{\sigma\left(  i\right)  }-x_{\sigma\left(  j\right)  }\right)
=\left(  -1\right)  ^{\sigma}\cdot\prod_{1\leq i<j\leq n}\left(  x_{i}%
-x_{j}\right)  $. This solves Exercise \ref{exe.perm.sign.pseudoexplicit}
\textbf{(a)}.

\textbf{(c)} Applying Exercise \ref{exe.perm.sign.pseudoexplicit} \textbf{(a)}
to $x_{i}=i$, we obtain $\prod_{1\leq i<j\leq n}\left(  \sigma\left(
i\right)  -\sigma\left(  j\right)  \right)  =\left(  -1\right)  ^{\sigma}%
\cdot\prod_{1\leq i<j\leq n}\left(  i-j\right)  $. We can divide both sides of
this equality by $\prod_{1\leq i<j\leq n}\left(  i-j\right)  $ (because
$\prod_{1\leq i<j\leq n}\left(  i-j\right)  $ is a product of nonzero
integers, and thus nonzero). As a result, we obtain $\dfrac{\prod_{1\leq
i<j\leq n}\left(  \sigma\left(  i\right)  -\sigma\left(  j\right)  \right)
}{\prod_{1\leq i<j\leq n}\left(  i-j\right)  }=\left(  -1\right)  ^{\sigma}$.
Thus,%
\[
\left(  -1\right)  ^{\sigma}=\dfrac{\prod_{1\leq i<j\leq n}\left(
\sigma\left(  i\right)  -\sigma\left(  j\right)  \right)  }{\prod_{1\leq
i<j\leq n}\left(  i-j\right)  }=\prod_{1\leq i<j\leq n}\dfrac{\sigma\left(
i\right)  -\sigma\left(  j\right)  }{i-j}.
\]
Thus, (\ref{eq.sign.pseudoexplicit}) is proven. This solves Exercise
\ref{exe.perm.sign.pseudoexplicit} \textbf{(c)}.

\textbf{(d)} See below.
\end{proof}
\end{verlong}

Exercise \ref{exe.perm.sign.pseudoexplicit} \textbf{(d)} asks us to give an
alternative solution to Exercise \ref{exe.ps2.2.5} \textbf{(b)}. Let us do
this now:

\begin{proof}
[Alternative solution to Exercise \ref{exe.ps2.2.5} \textbf{(b)}.]Let
$n\in\mathbb{N}$. Let $\sigma$ and $\tau$ be two permutations in $S_{n}$. We
need to show that $\ell\left(  \sigma\circ\tau\right)  \equiv\ell\left(
\sigma\right)  +\ell\left(  \tau\right)  \operatorname{mod}2$.

We are going to prove that $\left(  -1\right)  ^{\sigma\circ\tau}=\left(
-1\right)  ^{\sigma}\cdot\left(  -1\right)  ^{\tau}$ first.

Exercise \ref{exe.perm.sign.pseudoexplicit} \textbf{(a)} (applied to $x_{i}%
=i$) yields%
\[
\prod_{1\leq i<j\leq n}\left(  \sigma\left(  i\right)  -\sigma\left(
j\right)  \right)  =\left(  -1\right)  ^{\sigma}\cdot\prod_{1\leq i<j\leq
n}\left(  i-j\right)  .
\]
Exercise \ref{exe.perm.sign.pseudoexplicit} \textbf{(a)} (applied to
$\sigma\left(  i\right)  $ and $\tau$ instead of $x_{i}$ and $\sigma$) yields%
\begin{align}
\prod_{1\leq i<j\leq n}\left(  \sigma\left(  \tau\left(  i\right)  \right)
-\sigma\left(  \tau\left(  j\right)  \right)  \right)   &  =\left(  -1\right)
^{\tau}\cdot\underbrace{\prod_{1\leq i<j\leq n}\left(  \sigma\left(  i\right)
-\sigma\left(  j\right)  \right)  }_{=\left(  -1\right)  ^{\sigma}\cdot
\prod_{1\leq i<j\leq n}\left(  i-j\right)  }\nonumber\\
&  =\underbrace{\left(  -1\right)  ^{\tau}\cdot\left(  -1\right)  ^{\sigma}%
}_{=\left(  -1\right)  ^{\sigma}\cdot\left(  -1\right)  ^{\tau}}\cdot
\prod_{1\leq i<j\leq n}\left(  i-j\right) \nonumber\\
&  =\left(  -1\right)  ^{\sigma}\cdot\left(  -1\right)  ^{\tau}\cdot
\prod_{1\leq i<j\leq n}\left(  i-j\right)  .
\label{sol.perm.sign.pseudoexplicit.d.1}%
\end{align}
Exercise \ref{exe.perm.sign.pseudoexplicit} \textbf{(a)} (applied to
$\sigma\circ\tau$ instead of $\sigma$) yields%
\[
\prod_{1\leq i<j\leq n}\left(  \left(  \sigma\circ\tau\right)  \left(
i\right)  -\left(  \sigma\circ\tau\right)  \left(  j\right)  \right)  =\left(
-1\right)  ^{\sigma\circ\tau}\cdot\prod_{1\leq i<j\leq n}\left(  i-j\right)
.
\]
Thus,%
\begin{align}
\left(  -1\right)  ^{\sigma\circ\tau}\cdot\prod_{1\leq i<j\leq n}\left(
i-j\right)   &  =\prod_{1\leq i<j\leq n}\left(  \underbrace{\left(
\sigma\circ\tau\right)  \left(  i\right)  }_{=\sigma\left(  \tau\left(
i\right)  \right)  }-\underbrace{\left(  \sigma\circ\tau\right)  \left(
j\right)  }_{=\sigma\left(  \tau\left(  j\right)  \right)  }\right)
\nonumber\\
&  =\prod_{1\leq i<j\leq n}\left(  \sigma\left(  \tau\left(  i\right)
\right)  -\sigma\left(  \tau\left(  j\right)  \right)  \right) \nonumber\\
&  =\left(  -1\right)  ^{\sigma}\cdot\left(  -1\right)  ^{\tau}\cdot
\prod_{1\leq i<j\leq n}\left(  i-j\right)  \ \ \ \ \ \ \ \ \ \ \left(
\text{by (\ref{sol.perm.sign.pseudoexplicit.d.1})}\right)  .
\label{sol.perm.sign.pseudoexplicit.d.2}%
\end{align}

\begin{vershort}
But the integer $\prod_{1\leq i<j\leq n}\left(  i-j\right)  $ is nonzero
(since it is a product of the nonzero integers $i-j$). Hence, we can divide
both sides of the equality (\ref{sol.perm.sign.pseudoexplicit.d.2}) by
$\prod_{1\leq i<j\leq n}\left(  i-j\right)  $. We thus obtain $\left(
-1\right)  ^{\sigma\circ\tau}=\left(  -1\right)  ^{\sigma}\cdot\left(
-1\right)  ^{\tau}$.
\end{vershort}

\begin{verlong}
But the integer $\prod_{1\leq i<j\leq n}\left(  i-j\right)  $ is
nonzero\footnote{\textit{Proof.} If $i$ and $j$ are two integers satisfying
$1\leq i<j\leq n$, then the integer $i-j$ is negative (since $i<j$) and thus
nonzero. Hence, the product $\prod_{1\leq i<j\leq n}\left(  i-j\right)  $ is a
product of nonzero integers, and therefore itself nonzero (since a product of
nonzero integers is always nonzero). Qed.}. Hence, we can divide both sides of
the equality (\ref{sol.perm.sign.pseudoexplicit.d.2}) by $\prod_{1\leq i<j\leq
n}\left(  i-j\right)  $. We thus obtain $\left(  -1\right)  ^{\sigma\circ\tau
}=\left(  -1\right)  ^{\sigma}\cdot\left(  -1\right)  ^{\tau}$.
\end{verlong}

Now, the definition of $\left(  -1\right)  ^{\sigma}$ yields $\left(
-1\right)  ^{\sigma}=\left(  -1\right)  ^{\ell\left(  \sigma\right)  }$. Also,
the definition of $\left(  -1\right)  ^{\tau}$ yields $\left(  -1\right)
^{\tau}=\left(  -1\right)  ^{\ell\left(  \tau\right)  }$. Hence,
\[
\underbrace{\left(  -1\right)  ^{\sigma}}_{=\left(  -1\right)  ^{\ell\left(
\sigma\right)  }}\cdot\underbrace{\left(  -1\right)  ^{\tau}}_{=\left(
-1\right)  ^{\ell\left(  \tau\right)  }}=\left(  -1\right)  ^{\ell\left(
\sigma\right)  }\cdot\left(  -1\right)  ^{\ell\left(  \tau\right)  }=\left(
-1\right)  ^{\ell\left(  \sigma\right)  +\ell\left(  \tau\right)  }.
\]

Finally, the definition of $\left(  -1\right)  ^{\sigma\circ\tau}$ yields
$\left(  -1\right)  ^{\sigma\circ\tau}=\left(  -1\right)  ^{\ell\left(
\sigma\circ\tau\right)  }$. Thus,%
\[
\left(  -1\right)  ^{\ell\left(  \sigma\circ\tau\right)  }=\left(  -1\right)
^{\sigma\circ\tau}=\left(  -1\right)  ^{\sigma}\cdot\left(  -1\right)  ^{\tau
}=\left(  -1\right)  ^{\ell\left(  \sigma\right)  +\ell\left(  \tau\right)
}.
\]
But it is obvious that if two integers $u$ and $v$ satisfy $\left(  -1\right)
^{u}=\left(  -1\right)  ^{v}$, then $u\equiv v\operatorname{mod}%
2$\ \ \ \ \footnote{Indeed, Proposition \ref{prop.mod.parity} shows that the
two statements $\left(  u\equiv v\operatorname{mod}2\right)  $ and $\left(
\left(  -1\right)  ^{u}=\left(  -1\right)  ^{v}\right)  $ are equivalent.}.
Applying this to $u=\ell\left(  \sigma\circ\tau\right)  $ and $v=\ell\left(
\sigma\right)  +\ell\left(  \tau\right)  $, we obtain $\ell\left(  \sigma
\circ\tau\right)  \equiv\ell\left(  \sigma\right)  +\ell\left(  \tau\right)
\operatorname{mod}2$. Thus, Exercise \ref{exe.ps2.2.5} \textbf{(b)} is solved again.
\end{proof}

\subsection{Solution to Exercise \ref{exe.Ialbe}}

Exercise \ref{exe.Ialbe} is an example of a combinatorial fact that one can
easily convince oneself of (with some handwaving), but that is quite hard to
prove in a formal, bulletproof way. Thus, the solution given below is going to
be long, but we hope that the reader can avoid major parts of it by figuring
them out independently.

Before we start solving Exercise \ref{exe.Ialbe}, let us define some notations.

\begin{definition}
\label{def.sol.Ialbe.12n}For every $n\in\mathbb{N}$, we let $\left[  n\right]
$ denote the set $\left\{  1,2,\ldots,n\right\}  $.
\end{definition}

\begin{definition}
For every $n\in\mathbb{N}$ and every $n$-tuple $\mathbf{a}=\left(  a_{1}%
,a_{2},\ldots,a_{n}\right)  $ of integers, we define the following notations:

\begin{itemize}
\item An \textit{inversion} of $\mathbf{a}$ will mean a pair $\left(
i,j\right)  \in\left[  n\right]  ^{2}$ satisfying $i<j$ and $a_{i}>a_{j}$.

\item We denote by $\operatorname*{Inv}\left(  \mathbf{a}\right)  $ the set of
all inversions of $\mathbf{a}$. Thus, $\operatorname*{Inv}\left(
\mathbf{a}\right)  \subseteq\left[  n\right]  ^{2}$. More precisely,%
\begin{align}
\operatorname*{Inv}\left(  \mathbf{a}\right)   &  =\left(  \text{the set of
all inversions of }\mathbf{a}\right) \nonumber\\
&  =\left\{  \left(  i,j\right)  \in\left[  n\right]  ^{2}\ \mid\ i<j\text{
and }a_{i}>a_{j}\right\} \label{sol.Ialbe.Inv.1}\\
&  \ \ \ \ \ \ \ \ \ \ \left(
\begin{array}
[c]{c}%
\text{since the inversions of }\mathbf{a}\text{ are the pairs }\left(
i,j\right)  \in\left[  n\right]  ^{2}\\
\text{satisfying }i<j\text{ and }a_{i}>a_{j}\text{ (by the definition}\\
\text{of an \textquotedblleft inversion\textquotedblright)}%
\end{array}
\right) \nonumber\\
&  =\left\{  \left(  u,v\right)  \in\left[  n\right]  ^{2}\ \mid\ u<v\text{
and }a_{u}>a_{v}\right\} \label{sol.Ialbe.Inv.2}\\
&  \ \ \ \ \ \ \ \ \ \ \left(  \text{here, we renamed the index }\left(
i,j\right)  \text{ as }\left(  u,v\right)  \right)  .\nonumber
\end{align}

\item We denote by $\ell\left(  \mathbf{a}\right)  $ the number $\left\vert
\operatorname*{Inv}\left(  \mathbf{a}\right)  \right\vert $. (This is
well-defined because $\operatorname*{Inv}\left(  \mathbf{a}\right)  $ is
finite (since $\operatorname*{Inv}\left(  \mathbf{a}\right)  \subseteq\left[
n\right]  ^{2}$).) Thus,%
\begin{align}
\ell\left(  \mathbf{a}\right)   &  =\left\vert \underbrace{\operatorname*{Inv}%
\left(  \mathbf{a}\right)  }_{=\left(  \text{the set of all inversions of
}\mathbf{a}\right)  }\right\vert =\left\vert \left(  \text{the set of all
inversions of }\mathbf{a}\right)  \right\vert \nonumber\\
&  =\left(  \text{the number of all inversions of }\mathbf{a}\right)  .
\label{sol.Ialbe.l(a)}%
\end{align}

\end{itemize}
\end{definition}

The next, nearly trivial, lemma relies on the notion of the ``list of all
elements of $S$ in increasing order (with no repetitions)'', where $S$ is a
finite set of integers. This notion means exactly what it says (but see
Definition \ref{def.ind.inclist} for a rigorous definition).

\begin{lemma}
\label{lem.sol.Ialbe.Inv=Inv}Let $P$ be a finite set of integers. Let
$m=\left\vert P\right\vert $. Let $\sigma\in S_{m}$. Let $\left(  p_{1}%
,p_{2},\ldots,p_{m}\right)  $ be the list of all elements of $P$ in increasing
order (with no repetitions). Let $\operatorname*{Inv}\left(  \sigma\right)  $
denote the set of all inversions of $\sigma$.

\textbf{(a)} We have $\operatorname*{Inv}\left(  \sigma\right)
=\operatorname*{Inv}\left(  p_{\sigma\left(  1\right)  },p_{\sigma\left(
2\right)  },\ldots,p_{\sigma\left(  m\right)  }\right)  $.

\textbf{(b)} We have $\ell\left(  \sigma\right)  =\ell\left(  p_{\sigma\left(
1\right)  },p_{\sigma\left(  2\right)  },\ldots,p_{\sigma\left(  m\right)
}\right)  $.
\end{lemma}

\begin{vershort}
\begin{proof}
[Proof of Lemma \ref{lem.sol.Ialbe.Inv=Inv}.]The inversions of $\sigma$ are
the pairs $\left(  i,j\right)  $ of integers satisfying $1\leq i<j\leq m$ and
$\sigma\left(  i\right)  >\sigma\left(  j\right)  $. Thus,
$\operatorname*{Inv}\left(  \sigma\right)  $ (which is the set of all
inversions of $\sigma$) is the set of all such pairs $\left(  i,j\right)  $.
In other words,
\begin{align}
\operatorname*{Inv}\left(  \sigma\right)   &  =\left\{  \left(  i,j\right)
\in\mathbb{Z}^{2}\ \mid\ 1\leq i<j\leq m\text{ and }\sigma\left(  i\right)
>\sigma\left(  j\right)  \right\} \nonumber\\
&  =\left\{  \left(  i,j\right)  \in\left[  m\right]  ^{2}\ \mid\ i<j\text{
and }\sigma\left(  i\right)  >\sigma\left(  j\right)  \right\} \nonumber\\
&  \ \ \ \ \ \ \ \ \ \ \left(
\begin{array}
[c]{c}%
\text{because the pairs }\left(  i,j\right)  \in\mathbb{Z}^{2}\text{
satisfying }1\leq i<j\leq m\\
\text{are precisely the pairs }\left(  i,j\right)  \in\left[  m\right]
^{2}\text{ satisfying }i<j
\end{array}
\right) \nonumber\\
&  =\left\{  \left(  u,v\right)  \in\left[  m\right]  ^{2}\ \mid\ u<v\text{
and }\sigma\left(  u\right)  >\sigma\left(  v\right)  \right\}
\label{pf.lem.sol.Ialbe.Inv=Inv.short.lhs=}%
\end{align}
(here, we have renamed the index $\left(  i,j\right)  $ as $\left(
u,v\right)  $).

But $\left(  p_{1},p_{2},\ldots,p_{m}\right)  $ is the list of all elements of
$P$ in increasing order (with no repetitions). Thus, $\left(  p_{1}%
,p_{2},\ldots,p_{m}\right)  $ is a strictly increasing list. In other words,
$p_{1}<p_{2}<\cdots<p_{m}$. Hence, if $i$ and $j$ are two elements of $\left[
m\right]  $, then we have the following logical equivalence:%
\begin{equation}
\left(  p_{i}>p_{j}\right)  \Longleftrightarrow\left(  i>j\right)  .
\label{pf.lem.sol.Ialbe.Inv=Inv.short.1}%
\end{equation}

But (\ref{sol.Ialbe.Inv.2}) (applied to $m$, $\left(  p_{\sigma\left(
1\right)  },p_{\sigma\left(  2\right)  },\ldots,p_{\sigma\left(  m\right)
}\right)  $ and $p_{\sigma\left(  i\right)  }$ instead of $n$, $\mathbf{a}$
and $a_{i}$) yields%
\begin{align*}
&  \operatorname*{Inv}\left(  p_{\sigma\left(  1\right)  },p_{\sigma\left(
2\right)  },\ldots,p_{\sigma\left(  m\right)  }\right) \\
&  =\left\{  \left(  u,v\right)  \in\left[  m\right]  ^{2}\ \mid\ u<v\text{
and }\underbrace{p_{\sigma\left(  u\right)  }>p_{\sigma\left(  v\right)  }%
}_{\substack{\text{this is equivalent to }\sigma\left(  u\right)
>\sigma\left(  v\right)  \\\text{(by (\ref{pf.lem.sol.Ialbe.Inv=Inv.short.1}),
applied to }i=\sigma\left(  u\right)  \text{ and }j=\sigma\left(  v\right)
\text{)}}}\right\} \\
&  =\left\{  \left(  u,v\right)  \in\left[  m\right]  ^{2}\ \mid\ u<v\text{
and }\sigma\left(  u\right)  >\sigma\left(  v\right)  \right\}
=\operatorname*{Inv}\left(  \sigma\right)
\end{align*}
(by (\ref{pf.lem.sol.Ialbe.Inv=Inv.short.lhs=})). This proves Lemma
\ref{lem.sol.Ialbe.Inv=Inv} \textbf{(a)}.

\textbf{(b)} We know that $\ell\left(  \sigma\right)  $ is the number of all
inversions of $\sigma$ (by the definition of $\ell\left(  \sigma\right)  $).
In other words, $\ell\left(  \sigma\right)  $ is the size of the set of all
inversions of $\sigma$. In other words, $\ell\left(  \sigma\right)  $ is the
size of $\operatorname*{Inv}\left(  \sigma\right)  $ (since the set of all
inversions of $\sigma$ is $\operatorname*{Inv}\left(  \sigma\right)  $).
Hence,%
\begin{align*}
\ell\left(  \sigma\right)   &  =\left\vert \underbrace{\operatorname*{Inv}%
\left(  \sigma\right)  }_{\substack{=\operatorname*{Inv}\left(  p_{\sigma
\left(  1\right)  },p_{\sigma\left(  2\right)  },\ldots,p_{\sigma\left(
m\right)  }\right)  \\\text{(by Lemma \ref{lem.sol.Ialbe.Inv=Inv}
\textbf{(a)})}}}\right\vert =\left\vert \operatorname*{Inv}\left(
p_{\sigma\left(  1\right)  },p_{\sigma\left(  2\right)  },\ldots
,p_{\sigma\left(  m\right)  }\right)  \right\vert \\
&  =\ell\left(  p_{\sigma\left(  1\right)  },p_{\sigma\left(  2\right)
},\ldots,p_{\sigma\left(  m\right)  }\right)
\end{align*}
(since $\ell\left(  p_{\sigma\left(  1\right)  },p_{\sigma\left(  2\right)
},\ldots,p_{\sigma\left(  m\right)  }\right)  $ is defined as $\left\vert
\operatorname*{Inv}\left(  p_{\sigma\left(  1\right)  },p_{\sigma\left(
2\right)  },\ldots,p_{\sigma\left(  m\right)  }\right)  \right\vert $). This
proves Lemma \ref{lem.sol.Ialbe.Inv=Inv} \textbf{(b)}.
\end{proof}
\end{vershort}

\begin{verlong}
\begin{proof}
[Proof of Lemma \ref{lem.sol.Ialbe.Inv=Inv}.]We know that $\left(  p_{1}%
,p_{2},\ldots,p_{m}\right)  $ is the list of all elements of $P$ in increasing
order (with no repetitions). Thus, $\left(  p_{1},p_{2},\ldots,p_{m}\right)  $
is a strictly increasing list. In other words, we have $p_{1}<p_{2}%
<\cdots<p_{m}$. In other words, if $u$ and $v$ are two elements of $\left[
m\right]  $ such that $u<v$, then%
\begin{equation}
p_{u}<p_{v}. \label{pf.lem.sol.Ialbe.Inv=Inv.pupv}%
\end{equation}

\textbf{(a)} We have%
\begin{equation}
\operatorname*{Inv}\left(  p_{\sigma\left(  1\right)  },p_{\sigma\left(
2\right)  },\ldots,p_{\sigma\left(  m\right)  }\right)  =\left\{  \left(
u,v\right)  \in\left[  m\right]  ^{2}\ \mid\ u<v\text{ and }p_{\sigma\left(
u\right)  }>p_{\sigma\left(  v\right)  }\right\}
\label{pf.lem.sol.Ialbe.Inv=Inv.a.Inv}%
\end{equation}
(by (\ref{sol.Ialbe.Inv.2}), applied to $m$, $\left(  p_{\sigma\left(
1\right)  },p_{\sigma\left(  2\right)  },\ldots,p_{\sigma\left(  m\right)
}\right)  $ and $p_{\sigma\left(  i\right)  }$ instead of $n$, $\mathbf{a}$
and $a_{i}$).

Let $c\in\operatorname*{Inv}\left(  \sigma\right)  $. Thus, $c$ is an element
of $\operatorname*{Inv}\left(  \sigma\right)  $. In other words, $c$ is an
inversion of $\sigma$ (since $\operatorname*{Inv}\left(  \sigma\right)  $ is
the set of all inversions of $\sigma$). In other words, $c$ is a pair $\left(
i,j\right)  $ of integers satisfying $1\leq i<j\leq m$ and $\sigma\left(
i\right)  >\sigma\left(  j\right)  $ (by the definition of an
\textquotedblleft inversion of $\sigma$\textquotedblright).

From $i<j\leq m$, we obtain $i\leq m$. From $1\leq i\leq m$, we obtain
$i\in\left[  m\right]  $.

From $1\leq i<j$, we obtain $1\leq j$. From $1\leq j\leq m$, we obtain
$j\in\left[  m\right]  $.

From $i\in\left[  m\right]  $ and $j\in\left[  m\right]  $, we obtain $\left(
i,j\right)  \in\left[  m\right]  \times\left[  m\right]  =\left[  m\right]
^{2}$. From $\sigma\left(  i\right)  >\sigma\left(  j\right)  $, we obtain
$\sigma\left(  j\right)  <\sigma\left(  i\right)  $. Hence, $p_{\sigma\left(
j\right)  }<p_{\sigma\left(  i\right)  }$ (by
(\ref{pf.lem.sol.Ialbe.Inv=Inv.pupv}), applied to $u=\sigma\left(  j\right)  $
and $v=\sigma\left(  i\right)  $). In other words, $p_{\sigma\left(  i\right)
}>p_{\sigma\left(  j\right)  }$. Hence, $\left(  i,j\right)  $ is an element
of $\left[  m\right]  ^{2}$ and satisfies $i<j$ and $p_{\sigma\left(
i\right)  }>p_{\sigma\left(  j\right)  }$. In other words, $\left(
i,j\right)  $ is an element $\left(  u,v\right)  \in\left[  m\right]  ^{2}$
satisfying $u<v$ and $p_{\sigma\left(  u\right)  }>p_{\sigma\left(  v\right)
}$. Thus,%
\begin{align*}
\left(  i,j\right)   &  \in\left\{  \left(  u,v\right)  \in\left[  m\right]
^{2}\ \mid\ u<v\text{ and }p_{\sigma\left(  u\right)  }>p_{\sigma\left(
v\right)  }\right\} \\
&  =\operatorname*{Inv}\left(  p_{\sigma\left(  1\right)  },p_{\sigma\left(
2\right)  },\ldots,p_{\sigma\left(  m\right)  }\right)
\ \ \ \ \ \ \ \ \ \ \left(  \text{by (\ref{pf.lem.sol.Ialbe.Inv=Inv.a.Inv}%
)}\right)  .
\end{align*}
Hence, $c=\left(  i,j\right)  \in\operatorname*{Inv}\left(  p_{\sigma\left(
1\right)  },p_{\sigma\left(  2\right)  },\ldots,p_{\sigma\left(  m\right)
}\right)  $.

Now, forget that we fixed $c$. We thus have shown that every $c\in
\operatorname*{Inv}\left(  \sigma\right)  $ satisfies $c\in\operatorname*{Inv}%
\left(  p_{\sigma\left(  1\right)  },p_{\sigma\left(  2\right)  }%
,\ldots,p_{\sigma\left(  m\right)  }\right)  $. In other words,
\begin{equation}
\operatorname*{Inv}\left(  \sigma\right)  \subseteq\operatorname*{Inv}\left(
p_{\sigma\left(  1\right)  },p_{\sigma\left(  2\right)  },\ldots
,p_{\sigma\left(  m\right)  }\right)  .
\label{pf.lem.sol.Ialbe.Inv=Inv.a.dir1}%
\end{equation}

On the other hand, let $d\in\operatorname*{Inv}\left(  p_{\sigma\left(
1\right)  },p_{\sigma\left(  2\right)  },\ldots,p_{\sigma\left(  m\right)
}\right)  $. Thus,%
\begin{align*}
d  &  \in\operatorname*{Inv}\left(  p_{\sigma\left(  1\right)  }%
,p_{\sigma\left(  2\right)  },\ldots,p_{\sigma\left(  m\right)  }\right) \\
&  =\left\{  \left(  u,v\right)  \in\left[  m\right]  ^{2}\ \mid\ u<v\text{
and }p_{\sigma\left(  u\right)  }>p_{\sigma\left(  v\right)  }\right\}
\end{align*}
(by (\ref{pf.lem.sol.Ialbe.Inv=Inv.a.Inv})). In other words, $d=\left(
u,v\right)  $ for some $\left(  u,v\right)  \in\left[  m\right]  ^{2}$
satisfying $u<v$ and $p_{\sigma\left(  u\right)  }>p_{\sigma\left(  v\right)
}$. Consider this $\left(  u,v\right)  $.

From $\left(  u,v\right)  \in\left[  m\right]  ^{2}$, we obtain $u\in\left[
m\right]  $ and $v\in\left[  m\right]  $. From $u\in\left[  m\right]  $, we
obtain $1\leq u$. From $v\in\left[  m\right]  $, we obtain $v\leq m$. Thus,
$1\leq u<v\leq m$.

From $u<v$, we obtain $u\neq v$. If we had $\sigma\left(  u\right)
=\sigma\left(  v\right)  $, then we would have $\sigma^{-1}\left(
\underbrace{\sigma\left(  u\right)  }_{=\sigma\left(  v\right)  }\right)
=\sigma^{-1}\left(  \sigma\left(  v\right)  \right)  =v$, which would
contradict $\sigma^{-1}\left(  \sigma\left(  u\right)  \right)  =u\neq v$.
Thus, we cannot have $\sigma\left(  u\right)  =\sigma\left(  v\right)  $. In
other words, we have $\sigma\left(  u\right)  \neq\sigma\left(  v\right)  $.

Now, let us assume (for the sake of contradiction) that $\sigma\left(
u\right)  <\sigma\left(  v\right)  $. Then, $p_{\sigma\left(  u\right)
}<p_{\sigma\left(  v\right)  }$ (by (\ref{pf.lem.sol.Ialbe.Inv=Inv.pupv}),
applied to $\sigma\left(  u\right)  $ and $\sigma\left(  v\right)  $ instead
of $u$ and $v$). This contradicts $p_{\sigma\left(  u\right)  }>p_{\sigma
\left(  v\right)  }$. This contradiction proves that our assumption (that
$\sigma\left(  u\right)  <\sigma\left(  v\right)  $) was false. Hence, we
cannot have $\sigma\left(  u\right)  <\sigma\left(  v\right)  $. We thus have
$\sigma\left(  u\right)  \geq\sigma\left(  v\right)  $. Combining this with
$\sigma\left(  u\right)  \neq\sigma\left(  v\right)  $, we obtain
$\sigma\left(  u\right)  >\sigma\left(  v\right)  $.

So we know that $\left(  u,v\right)  $ is a pair of integers satisfying $1\leq
u<v\leq m$ and $\sigma\left(  u\right)  >\sigma\left(  v\right)  $. In other
words, $\left(  u,v\right)  $ is a pair $\left(  i,j\right)  $ of integers
satisfying $1\leq i<j\leq m$ and $\sigma\left(  i\right)  >\sigma\left(
j\right)  $. In other words, $\left(  u,v\right)  $ is an inversion of
$\sigma$ (by the definition of an \textquotedblleft inversion of $\sigma
$\textquotedblright). In other words, $\left(  u,v\right)  $ is an element of
$\operatorname*{Inv}\left(  \sigma\right)  $ (since $\operatorname*{Inv}%
\left(  \sigma\right)  $ is the set of all inversions of $\sigma$). In other
words, $\left(  u,v\right)  \in\operatorname*{Inv}\left(  \sigma\right)  $.
Thus, $d=\left(  u,v\right)  \in\operatorname*{Inv}\left(  \sigma\right)  $.

Now, forget that we fixed $d$. We thus have shown that every \newline%
$d\in\operatorname*{Inv}\left(  p_{\sigma\left(  1\right)  },p_{\sigma\left(
2\right)  },\ldots,p_{\sigma\left(  m\right)  }\right)  $ satisfies
$d\in\operatorname*{Inv}\left(  \sigma\right)  $. In other words,%
\[
\operatorname*{Inv}\left(  p_{\sigma\left(  1\right)  },p_{\sigma\left(
2\right)  },\ldots,p_{\sigma\left(  m\right)  }\right)  \subseteq
\operatorname*{Inv}\left(  \sigma\right)  .
\]
Combining this with (\ref{pf.lem.sol.Ialbe.Inv=Inv.a.dir1}), we obtain
$\operatorname*{Inv}\left(  \sigma\right)  =\operatorname*{Inv}\left(
p_{\sigma\left(  1\right)  },p_{\sigma\left(  2\right)  },\ldots
,p_{\sigma\left(  m\right)  }\right)  $. This proves Lemma
\ref{lem.sol.Ialbe.Inv=Inv} \textbf{(a)}.

\textbf{(b)} Recall that $\ell\left(  \sigma\right)  $ is the number of
inversions of $\sigma$ (indeed, this is how $\ell\left(  \sigma\right)  $ was
defined). Thus,%
\begin{align*}
\ell\left(  \sigma\right)   &  =\left(  \text{the number of inversions of
}\sigma\right) \\
&  =\left\vert \underbrace{\left(  \text{the set of all inversions of }%
\sigma\right)  }_{\substack{=\operatorname*{Inv}\left(  \sigma\right)
\\\text{(since }\operatorname*{Inv}\left(  \sigma\right)  \text{ is the set of
all inversions of }\sigma\text{)}}}\right\vert \\
&  =\left\vert \underbrace{\operatorname*{Inv}\left(  \sigma\right)
}_{\substack{=\operatorname*{Inv}\left(  p_{\sigma\left(  1\right)
},p_{\sigma\left(  2\right)  },\ldots,p_{\sigma\left(  m\right)  }\right)
\\\text{(by Lemma \ref{lem.sol.Ialbe.Inv=Inv} \textbf{(a)})}}}\right\vert
=\left\vert \operatorname*{Inv}\left(  p_{\sigma\left(  1\right)  }%
,p_{\sigma\left(  2\right)  },\ldots,p_{\sigma\left(  m\right)  }\right)
\right\vert .
\end{align*}
Comparing this with%
\begin{align*}
\ell\left(  p_{\sigma\left(  1\right)  },p_{\sigma\left(  2\right)  }%
,\ldots,p_{\sigma\left(  m\right)  }\right)   &  =\left\vert
\operatorname*{Inv}\left(  p_{\sigma\left(  1\right)  },p_{\sigma\left(
2\right)  },\ldots,p_{\sigma\left(  m\right)  }\right)  \right\vert \\
&  \ \ \ \ \ \ \ \ \ \ \left(  \text{by the definition of }\ell\left(
p_{\sigma\left(  1\right)  },p_{\sigma\left(  2\right)  },\ldots
,p_{\sigma\left(  m\right)  }\right)  \right)  ,
\end{align*}
we obtain $\ell\left(  \sigma\right)  =\ell\left(  p_{\sigma\left(  1\right)
},p_{\sigma\left(  2\right)  },\ldots,p_{\sigma\left(  m\right)  }\right)  $.
This proves Lemma \ref{lem.sol.Ialbe.Inv=Inv} \textbf{(b)}.
\end{proof}
\end{verlong}

In the following, we shall use the Iverson bracket notation introduced in
Definition \ref{def.iverson}.

Now, let us show some more lemmas:

\begin{lemma}
\label{lem.sol.Ialbe.iverson}Let $n\in\mathbb{N}$. Let $\mathbf{a}=\left(
a_{1},a_{2},\ldots,a_{n}\right)  $ be an $n$-tuple of integers. Then,%
\[
\ell\left(  \mathbf{a}\right)  =\sum_{i=1}^{n}\sum_{j=1}^{n}\left[
i<j\right]  \left[  a_{i}>a_{j}\right]  .
\]

\end{lemma}

\begin{vershort}
\begin{proof}
[Proof of Lemma \ref{lem.sol.Ialbe.iverson}.]If $\mathcal{A}$ and
$\mathcal{B}$ are two logical statements, then%
\begin{equation}
\left[  \mathcal{A}\right]  \left[  \mathcal{B}\right]  =\left[
\mathcal{A}\text{ and }\mathcal{B}\right]
\label{pf.lem.sol.Ialbe.iverson.short.AandB}%
\end{equation}
(by Exercise \ref{exe.iverson-prop} \textbf{(b)}).

We have%
\begin{align*}
&  \underbrace{\sum_{i=1}^{n}\sum_{j=1}^{n}}_{=\sum_{i\in\left[  n\right]
}\sum_{j\in\left[  n\right]  }=\sum_{\left(  i,j\right)  \in\left[  n\right]
^{2}}}\underbrace{\left[  i<j\right]  \left[  a_{i}>a_{j}\right]
}_{\substack{=\left[  i<j\text{ and }a_{i}>a_{j}\right]  \\\text{(by
(\ref{pf.lem.sol.Ialbe.iverson.short.AandB}))}}}\\
&  =\sum_{\left(  i,j\right)  \in\left[  n\right]  ^{2}}\left[  i<j\text{ and
}a_{i}>a_{j}\right] \\
&  =\sum_{\substack{\left(  i,j\right)  \in\left[  n\right]  ^{2};\\i<j\text{
and }a_{i}>a_{j}}}\underbrace{\left[  i<j\text{ and }a_{i}>a_{j}\right]
}_{\substack{=1\\\text{(since }i<j\text{ and }a_{i}>a_{j}\text{)}}%
}+\sum_{\substack{\left(  i,j\right)  \in\left[  n\right]  ^{2};\\\text{not
}\left(  i<j\text{ and }a_{i}>a_{j}\right)  }}\underbrace{\left[  i<j\text{
and }a_{i}>a_{j}\right]  }_{\substack{=0\\\text{(since not }\left(  i<j\text{
and }a_{i}>a_{j}\right)  \text{)}}}\\
&  =\sum_{\substack{\left(  i,j\right)  \in\left[  n\right]  ^{2};\\i<j\text{
and }a_{i}>a_{j}}}1+\underbrace{\sum_{\substack{\left(  i,j\right)  \in\left[
n\right]  ^{2};\\\text{not }\left(  i<j\text{ and }a_{i}>a_{j}\right)  }%
}0}_{=0}=\sum_{\substack{\left(  i,j\right)  \in\left[  n\right]
^{2};\\i<j\text{ and }a_{i}>a_{j}}}1\\
&  =\left\vert \underbrace{\left\{  \left(  i,j\right)  \in\left[  n\right]
^{2}\ \mid\ i<j\text{ and }a_{i}>a_{j}\right\}  }%
_{\substack{=\operatorname*{Inv}\left(  \mathbf{a}\right)  \\\text{(by
(\ref{sol.Ialbe.Inv.1}))}}}\right\vert \cdot1=\left\vert \operatorname*{Inv}%
\left(  \mathbf{a}\right)  \right\vert \cdot1\\
&  =\left\vert \operatorname*{Inv}\left(  \mathbf{a}\right)  \right\vert
=\ell\left(  \mathbf{a}\right)  \ \ \ \ \ \ \ \ \ \ \left(  \text{since }%
\ell\left(  \mathbf{a}\right)  \text{ was defined to be }\left\vert
\operatorname*{Inv}\left(  \mathbf{a}\right)  \right\vert \right)  .
\end{align*}
This proves Lemma \ref{lem.sol.Ialbe.iverson}.
\end{proof}
\end{vershort}

\begin{verlong}
\begin{proof}
[Proof of Lemma \ref{lem.sol.Ialbe.iverson}.]We have%
\begin{align*}
\sum_{c\in\operatorname*{Inv}\left(  \mathbf{a}\right)  }1  &  =\sum
_{c\in\left\{  \left(  u,v\right)  \in\left[  n\right]  ^{2}\ \mid\ u<v\text{
and }a_{u}>a_{v}\right\}  }1\\
&  \ \ \ \ \ \ \ \ \ \ \left(  \text{since }\operatorname*{Inv}\left(
\mathbf{a}\right)  =\left\{  \left(  u,v\right)  \in\left[  n\right]
^{2}\ \mid\ u<v\text{ and }a_{u}>a_{v}\right\}  \text{ (by
(\ref{sol.Ialbe.Inv.2}))}\right) \\
&  =\underbrace{\sum_{\left(  i,j\right)  \in\left\{  \left(  u,v\right)
\in\left[  n\right]  ^{2}\ \mid\ u<v\text{ and }a_{u}>a_{v}\right\}  }}%
_{=\sum_{\substack{\left(  i,j\right)  \in\left[  n\right]  ^{2};\\i<j\text{
and }a_{i}>a_{j}}}}1\\
&  \ \ \ \ \ \ \ \ \ \ \left(  \text{here, we have substituted }\left(
i,j\right)  \text{ for }c\text{ in the sum}\right) \\
&  =\sum_{\substack{\left(  i,j\right)  \in\left[  n\right]  ^{2};\\i<j\text{
and }a_{i}>a_{j}}}1.
\end{align*}
Comparing this with%
\begin{align*}
&  \sum_{\substack{\left(  i,j\right)  \in\left[  n\right]  ^{2}%
;\\i<j}}\left[  a_{i}>a_{j}\right] \\
&  =\sum_{\substack{\left(  i,j\right)  \in\left[  n\right]  ^{2};\\i<j\text{
and }a_{i}>a_{j}}}\underbrace{\left[  a_{i}>a_{j}\right]  }%
_{\substack{=1\\\text{(since }a_{i}>a_{j}\text{)}}}+\sum_{\substack{\left(
i,j\right)  \in\left[  n\right]  ^{2};\\i<j\text{ and not }a_{i}>a_{j}%
}}\underbrace{\left[  a_{i}>a_{j}\right]  }_{\substack{=0\\\text{(since we
don't have }a_{i}>a_{j}\text{)}}}\\
&  \ \ \ \ \ \ \ \ \ \ \left(
\begin{array}
[c]{c}%
\text{since every }\left(  i,j\right)  \in\left[  n\right]  ^{2}\text{
satisfies either }a_{i}>a_{j}\text{ or}\\
\left(  \text{not }a_{i}>a_{j}\right)  \text{ (but not both)}%
\end{array}
\right) \\
&  =\sum_{\substack{\left(  i,j\right)  \in\left[  n\right]  ^{2};\\i<j\text{
and }a_{i}>a_{j}}}1+\underbrace{\sum_{\substack{\left(  i,j\right)  \in\left[
n\right]  ^{2};\\i<j\text{ and not }a_{i}>a_{j}}}0}_{=0}=\sum
_{\substack{\left(  i,j\right)  \in\left[  n\right]  ^{2};\\i<j\text{ and
}a_{i}>a_{j}}}1,
\end{align*}
we obtain%
\[
\sum_{\substack{\left(  i,j\right)  \in\left[  n\right]  ^{2};\\i<j}}\left[
a_{i}>a_{j}\right]  =\sum_{c\in\operatorname*{Inv}\left(  \mathbf{a}\right)
}1=\left\vert \operatorname*{Inv}\left(  \mathbf{a}\right)  \right\vert
\cdot1=\left\vert \operatorname*{Inv}\left(  \mathbf{a}\right)  \right\vert .
\]
Comparing this with%
\[
\ell\left(  \mathbf{a}\right)  =\left\vert \operatorname*{Inv}\left(
\mathbf{a}\right)  \right\vert \ \ \ \ \ \ \ \ \ \ \left(  \text{by the
definition of }\ell\left(  \mathbf{a}\right)  \right)  ,
\]
we obtain%
\[
\ell\left(  \mathbf{a}\right)  =\sum_{\substack{\left(  i,j\right)  \in\left[
n\right]  ^{2};\\i<j}}\left[  a_{i}>a_{j}\right]  .
\]
Comparing this with%
\begin{align*}
&  \sum_{\left(  i,j\right)  \in\left[  n\right]  ^{2}}\left[  i<j\right]
\left[  a_{i}>a_{j}\right] \\
&  =\sum_{\substack{\left(  i,j\right)  \in\left[  n\right]  ^{2}%
;\\i<j}}\underbrace{\left[  i<j\right]  }_{\substack{=1\\\text{(since
}i<j\text{)}}}\left[  a_{i}>a_{j}\right]  +\sum_{\substack{\left(  i,j\right)
\in\left[  n\right]  ^{2};\\\text{not }i<j}}\underbrace{\left[  i<j\right]
}_{\substack{=0\\\text{(since we don't have }i<j\text{)}}}\left[  a_{i}%
>a_{j}\right] \\
&  \ \ \ \ \ \ \ \ \ \ \left(
\begin{array}
[c]{c}%
\text{since every }\left(  i,j\right)  \in\left[  n\right]  ^{2}\text{
satisfies either }i<j\text{ or}\\
\left(  \text{not }i<j\right)  \text{ (but not both)}%
\end{array}
\right) \\
&  =\sum_{\substack{\left(  i,j\right)  \in\left[  n\right]  ^{2}%
;\\i<j}}\left[  a_{i}>a_{j}\right]  +\underbrace{\sum_{\substack{\left(
i,j\right)  \in\left[  n\right]  ^{2};\\\text{not }i<j}}0\left[  a_{i}%
>a_{j}\right]  }_{=0}=\sum_{\substack{\left(  i,j\right)  \in\left[  n\right]
^{2};\\i<j}}\left[  a_{i}>a_{j}\right]  ,
\end{align*}
we obtain%
\begin{align*}
\ell\left(  \mathbf{a}\right)   &  =\underbrace{\sum_{\left(  i,j\right)
\in\left[  n\right]  ^{2}}}_{=\sum_{i\in\left[  n\right]  }\sum_{j\in\left[
n\right]  }}\left[  i<j\right]  \left[  a_{i}>a_{j}\right] \\
&  =\underbrace{\sum_{i\in\left[  n\right]  }}_{\substack{=\sum_{i\in\left\{
1,2,\ldots,n\right\}  }\\\text{(since }\left[  n\right]  =\left\{
1,2,\ldots,n\right\}  \text{)}}}\underbrace{\sum_{j\in\left[  n\right]  }%
}_{\substack{=\sum_{j\in\left\{  1,2,\ldots,n\right\}  }\\\text{(since
}\left[  n\right]  =\left\{  1,2,\ldots,n\right\}  \text{)}}}\left[
i<j\right]  \left[  a_{i}>a_{j}\right] \\
&  =\underbrace{\sum_{i\in\left\{  1,2,\ldots,n\right\}  }}_{=\sum_{i=1}^{n}%
}\underbrace{\sum_{j\in\left\{  1,2,\ldots,n\right\}  }}_{=\sum_{j=1}^{n}%
}\left[  i<j\right]  \left[  a_{i}>a_{j}\right]  =\sum_{i=1}^{n}\sum_{j=1}%
^{n}\left[  i<j\right]  \left[  a_{i}>a_{j}\right]  .
\end{align*}
This proves Lemma \ref{lem.sol.Ialbe.iverson}.
\end{proof}
\end{verlong}

\begin{lemma}
\label{lem.sol.Ialbe.ab}Let $n\in\mathbb{N}$ and $m\in\mathbb{N}$. Let
$\mathbf{a}=\left(  a_{1},a_{2},\ldots,a_{n}\right)  $ be an $n$-tuple of
integers. Let $\mathbf{b}=\left(  b_{1},b_{2},\ldots,b_{m}\right)  $ be an
$m$-tuple of integers. Let $\mathbf{c}$ be the $\left(  n+m\right)  $-tuple
$\left(  a_{1},a_{2},\ldots,a_{n},b_{1},b_{2},\ldots,b_{m}\right)  $ of
integers. Then,%
\[
\ell\left(  \mathbf{c}\right)  =\ell\left(  \mathbf{a}\right)  +\ell\left(
\mathbf{b}\right)  +\sum_{\left(  i,j\right)  \in\left[  n\right]
\times\left[  m\right]  }\left[  a_{i}>b_{j}\right]  .
\]

\end{lemma}

\begin{vershort}
\begin{proof}
[Proof of Lemma \ref{lem.sol.Ialbe.ab}.]Lemma \ref{lem.sol.Ialbe.iverson}
yields%
\begin{equation}
\ell\left(  \mathbf{a}\right)  =\sum_{i=1}^{n}\sum_{j=1}^{n}\left[
i<j\right]  \left[  a_{i}>a_{j}\right]  .
\label{pf.lem.sol.Ialbe.ab.short.la=}%
\end{equation}
Also, Lemma \ref{lem.sol.Ialbe.iverson} (applied to $m$, $\mathbf{b}$ and
$b_{i}$ instead of $n$, $\mathbf{a}$ and $a_{i}$) yields%
\begin{equation}
\ell\left(  \mathbf{b}\right)  =\sum_{i=1}^{m}\sum_{j=1}^{m}\left[
i<j\right]  \left[  b_{i}>b_{j}\right]  .
\label{pf.lem.sol.Ialbe.ab.short.lb=}%
\end{equation}

Write the $\left(  n+m\right)  $-tuple $\mathbf{c}$ in the form $\mathbf{c}%
=\left(  c_{1},c_{2},\ldots,c_{n+m}\right)  $. Thus,%
\[
\left(  c_{1},c_{2},\ldots,c_{n+m}\right)  =\mathbf{c}=\left(  a_{1}%
,a_{2},\ldots,a_{n},b_{1},b_{2},\ldots,b_{m}\right)
\]
(by the definition of $\mathbf{c}$). In other words,%
\begin{equation}
\left(  c_{i}=a_{i}\ \ \ \ \ \ \ \ \ \ \text{for every }i\in\left\{
1,2,\ldots,n\right\}  \right)  \label{pf.lem.sol.Ialbe.ab.short.ci.a}%
\end{equation}
and%
\begin{equation}
\left(  c_{i}=b_{i-n}\ \ \ \ \ \ \ \ \ \ \text{for every }i\in\left\{
n+1,n+2,\ldots,n+m\right\}  \right)  . \label{pf.lem.sol.Ialbe.ab.short.ci.b}%
\end{equation}

But Lemma \ref{lem.sol.Ialbe.iverson} (applied to $n+m$, $\mathbf{c}$ and
$c_{i}$ instead of $n$, $\mathbf{a}$ and $a_{i}$) yields%
\begin{align}
\ell\left(  \mathbf{c}\right)   &  =\sum_{i=1}^{n+m}\sum_{j=1}^{n+m}\left[
i<j\right]  \left[  c_{i}>c_{j}\right] \nonumber\\
&  =\sum_{i=1}^{n}\sum_{j=1}^{n+m}\left[  i<j\right]  \left[  c_{i}%
>c_{j}\right]  +\sum_{i=n+1}^{n+m}\sum_{j=1}^{n+m}\left[  i<j\right]  \left[
c_{i}>c_{j}\right]  \label{pf.lem.sol.Ialbe.ab.short.1}%
\end{align}
(since $0\leq n\leq n+m$).

But every $i\in\left\{  1,2,\ldots,n\right\}  $ satisfies%
\begin{equation}
\sum_{j=1}^{n+m}\left[  i<j\right]  \left[  c_{i}>c_{j}\right]  =\sum
_{j=1}^{n}\left[  i<j\right]  \left[  a_{i}>a_{j}\right]  +\sum_{j=1}%
^{m}\left[  a_{i}>b_{j}\right]  \label{pf.lem.sol.Ialbe.ab.short.3a}%
\end{equation}
\footnote{\textit{Proof of (\ref{pf.lem.sol.Ialbe.ab.short.3a}):} Let
$i\in\left\{  1,2,\ldots,n\right\}  $. Thus, $c_{i}=a_{i}$ (by
(\ref{pf.lem.sol.Ialbe.ab.short.ci.a})).
\par
For every $j\in\left\{  1,2,\ldots,n\right\}  $, we have%
\begin{equation}
c_{j}=a_{j} \label{pf.lem.sol.Ialbe.ab.short.3a.pf.1}%
\end{equation}
(by (\ref{pf.lem.sol.Ialbe.ab.short.ci.a}), applied to $j$ instead of $i$).
\par
For every $j\in\left\{  n+1,n+2,\ldots,n+m\right\}  $, we have%
\begin{equation}
c_{j}=b_{j-n} \label{pf.lem.sol.Ialbe.ab.short.3a.pf.2}%
\end{equation}
(by (\ref{pf.lem.sol.Ialbe.ab.short.ci.b}), applied to $j$ instead of $i$).
\par
We have $i\leq n$ (since $i\in\left\{  1,2,\ldots,n\right\}  $). For every
$j\in\left\{  n+1,n+2,\ldots,n+m\right\}  $, we have $j\geq n+1>n\geq i$
(since $i\leq n$), thus $i<j$, therefore%
\begin{equation}
\left[  i<j\right]  =1. \label{pf.lem.sol.Ialbe.ab.short.3a.pf.3}%
\end{equation}
\par
Recall that $0\leq n\leq n+m$. Hence,%
\begin{align*}
\sum_{j=1}^{n+m}\left[  i<j\right]  \left[  c_{i}>c_{j}\right]   &
=\sum_{j=1}^{n}\left[  i<j\right]  \left[  \underbrace{c_{i}}_{=a_{i}%
}>\underbrace{c_{j}}_{\substack{=a_{j}\\\text{(by
(\ref{pf.lem.sol.Ialbe.ab.short.3a.pf.1}))}}}\right]  +\sum_{j=n+1}%
^{n+m}\underbrace{\left[  i<j\right]  }_{\substack{=1\\\text{(by
(\ref{pf.lem.sol.Ialbe.ab.short.3a.pf.3}))}}}\left[  \underbrace{c_{i}%
}_{=a_{i}}>\underbrace{c_{j}}_{\substack{=b_{j-n}\\\text{(by
(\ref{pf.lem.sol.Ialbe.ab.short.3a.pf.2}))}}}\right] \\
&  =\sum_{j=1}^{n}\left[  i<j\right]  \left[  a_{i}>a_{j}\right]
+\sum_{j=n+1}^{n+m}\left[  a_{i}>b_{j-n}\right] \\
&  =\sum_{j=1}^{n}\left[  i<j\right]  \left[  a_{i}>a_{j}\right]  +\sum
_{j=1}^{m}\left[  a_{i}>b_{j}\right] \\
&  \ \ \ \ \ \ \ \ \ \ \left(  \text{here, we have substituted }j\text{ for
}j-n\text{ in the second sum}\right)  .
\end{align*}
This proves (\ref{pf.lem.sol.Ialbe.ab.short.3a}).}.

Also, every $i\in\left\{  n+1,n+2,\ldots,n+m\right\}  $ satisfies%
\begin{equation}
\sum_{j=1}^{n+m}\left[  i<j\right]  \left[  c_{i}>c_{j}\right]  =\sum
_{j=1}^{m}\left[  i-n<j\right]  \left[  b_{i-n}>b_{j}\right]
\label{pf.lem.sol.Ialbe.ab.short.3b}%
\end{equation}
\footnote{\textit{Proof of (\ref{pf.lem.sol.Ialbe.ab.short.3b}):} Let
$i\in\left\{  n+1,n+2,\ldots,n+m\right\}  $. Thus, $c_{i}=b_{i-n}$ (by
(\ref{pf.lem.sol.Ialbe.ab.short.ci.b})).
\par
For every $j\in\left\{  n+1,n+2,\ldots,n+m\right\}  $, we have%
\begin{equation}
c_{j}=b_{j-n} \label{pf.lem.sol.Ialbe.ab.short.3b.pf.2}%
\end{equation}
(by (\ref{pf.lem.sol.Ialbe.ab.short.ci.b}), applied to $j$ instead of $i$).
\par
We have $i\geq n+1$ (since $i\in\left\{  n+1,n+2,\ldots,n+m\right\}  $) and
thus $i\geq n+1>n$, so that $n<i$. For every $j\in\left\{  1,2,\ldots
,n\right\}  $, we have $j\leq n\leq i$ (since $i\geq n$) and thus $i\geq j$.
Hence, for every $j\in\left\{  1,2,\ldots,n\right\}  $, we don't have $i<j$.
Thus, for every $j\in\left\{  1,2,\ldots,n\right\}  $, we have%
\begin{equation}
\left[  i<j\right]  =0. \label{pf.lem.sol.Ialbe.ab.short.3b.pf.3}%
\end{equation}
\par
Recall that $0\leq n\leq n+m$. Hence,%
\begin{align*}
\sum_{j=1}^{n+m}\left[  i<j\right]  \left[  c_{i}>c_{j}\right]   &
=\sum_{j=1}^{n}\underbrace{\left[  i<j\right]  }_{\substack{=0\\\text{(by
(\ref{pf.lem.sol.Ialbe.ab.short.3b.pf.3}))}}}\left[  c_{i}>c_{j}\right]
+\sum_{j=n+1}^{n+m}\underbrace{\left[  i<j\right]  }_{\substack{=\left[
i-n<j-n\right]  \\\text{(since }i<j\text{ is equivalent to}\\i-n<j-n\text{)}%
}}\left[  \underbrace{c_{i}}_{=b_{i-n}}>\underbrace{c_{j}}_{\substack{=b_{j-n}%
\\\text{(by (\ref{pf.lem.sol.Ialbe.ab.short.3b.pf.2}))}}}\right] \\
&  =\underbrace{\sum_{j=1}^{n}0\left[  c_{i}>c_{j}\right]  }_{=0}+\sum
_{j=n+1}^{n+m}\left[  i-n<j-n\right]  \left[  b_{i-n}>b_{j-n}\right] \\
&  =\sum_{j=n+1}^{n+m}\left[  i-n<j-n\right]  \left[  b_{i-n}>b_{j-n}\right]
\\
&  =\sum_{j=1}^{m}\left[  i-n<j\right]  \left[  b_{i-n}>b_{j}\right] \\
&  \ \ \ \ \ \ \ \ \ \ \left(  \text{here, we have substituted }j\text{ for
}j-n\text{ in the sum}\right)  .
\end{align*}
This proves (\ref{pf.lem.sol.Ialbe.ab.short.3b}).}.

Now, (\ref{pf.lem.sol.Ialbe.ab.short.1}) becomes%
\begin{align*}
&  \ell\left(  \mathbf{c}\right) \\
&  =\sum_{i=1}^{n}\underbrace{\sum_{j=1}^{n+m}\left[  i<j\right]  \left[
c_{i}>c_{j}\right]  }_{\substack{=\sum_{j=1}^{n}\left[  i<j\right]  \left[
a_{i}>a_{j}\right]  +\sum_{j=1}^{m}\left[  a_{i}>b_{j}\right]  \\\text{(by
(\ref{pf.lem.sol.Ialbe.ab.short.3a}))}}}+\sum_{i=n+1}^{n+m}\underbrace{\sum
_{j=1}^{n+m}\left[  i<j\right]  \left[  c_{i}>c_{j}\right]  }_{\substack{=\sum
_{j=1}^{m}\left[  i-n<j\right]  \left[  b_{i-n}>b_{j}\right]  \\\text{(by
(\ref{pf.lem.sol.Ialbe.ab.short.3b}))}}}\\
&  =\underbrace{\sum_{i=1}^{n}\left(  \sum_{j=1}^{n}\left[  i<j\right]
\left[  a_{i}>a_{j}\right]  +\sum_{j=1}^{m}\left[  a_{i}>b_{j}\right]
\right)  }_{=\sum_{i=1}^{n}\sum_{j=1}^{n}\left[  i<j\right]  \left[
a_{i}>a_{j}\right]  +\sum_{i=1}^{n}\sum_{j=1}^{m}\left[  a_{i}>b_{j}\right]
}+\underbrace{\sum_{i=n+1}^{n+m}\sum_{j=1}^{m}\left[  i-n<j\right]  \left[
b_{i-n}>b_{j}\right]  }_{\substack{=\sum_{i=1}^{m}\sum_{j=1}^{m}\left[
i<j\right]  \left[  b_{i}>b_{j}\right]  \\\text{(here, we have substituted
}i\text{ for }i-n\\\text{in the outer sum)}}}\\
&  =\underbrace{\sum_{i=1}^{n}\sum_{j=1}^{n}\left[  i<j\right]  \left[
a_{i}>a_{j}\right]  }_{\substack{=\ell\left(  \mathbf{a}\right)  \\\text{(by
(\ref{pf.lem.sol.Ialbe.ab.short.la=}))}}}+\underbrace{\sum_{i=1}^{n}%
}_{\substack{=\sum_{i\in\left[  n\right]  }}}\underbrace{\sum_{j=1}^{m}%
}_{\substack{=\sum_{j\in\left[  m\right]  }}}\left[  a_{i}>b_{j}\right]
+\underbrace{\sum_{i=1}^{m}\sum_{j=1}^{m}\left[  i<j\right]  \left[
b_{i}>b_{j}\right]  }_{\substack{=\ell\left(  \mathbf{b}\right)  \\\text{(by
(\ref{pf.lem.sol.Ialbe.ab.short.lb=}))}}}\\
&  =\ell\left(  \mathbf{a}\right)  +\underbrace{\sum_{i\in\left[  n\right]
}\sum_{j\in\left[  m\right]  }}_{=\sum_{\left(  i,j\right)  \in\left[
n\right]  \times\left[  m\right]  }}\left[  a_{i}>b_{j}\right]  +\ell\left(
\mathbf{b}\right) \\
&  =\ell\left(  \mathbf{a}\right)  +\sum_{\left(  i,j\right)  \in\left[
n\right]  \times\left[  m\right]  }\left[  a_{i}>b_{j}\right]  +\ell\left(
\mathbf{b}\right)  =\ell\left(  \mathbf{a}\right)  +\ell\left(  \mathbf{b}%
\right)  +\sum_{\left(  i,j\right)  \in\left[  n\right]  \times\left[
m\right]  }\left[  a_{i}>b_{j}\right]  .
\end{align*}
This proves Lemma \ref{lem.sol.Ialbe.ab}.
\end{proof}
\end{vershort}

\begin{verlong}
\begin{proof}
[Proof of Lemma \ref{lem.sol.Ialbe.ab}.]We have $n+\underbrace{m}_{\geq0}\geq
n$, so that $n\leq n+m$. Also, $n\geq0$, so that $0\leq n$.

Lemma \ref{lem.sol.Ialbe.iverson} yields%
\begin{equation}
\ell\left(  \mathbf{a}\right)  =\sum_{i=1}^{n}\sum_{j=1}^{n}\left[
i<j\right]  \left[  a_{i}>a_{j}\right]  . \label{pf.lem.sol.Ialbe.ab.la=}%
\end{equation}
Also, Lemma \ref{lem.sol.Ialbe.iverson} (applied to $m$, $\mathbf{b}$ and
$b_{i}$ instead of $n$, $\mathbf{a}$ and $a_{i}$) yields%
\begin{equation}
\ell\left(  \mathbf{b}\right)  =\sum_{i=1}^{m}\sum_{j=1}^{m}\left[
i<j\right]  \left[  b_{i}>b_{j}\right]  . \label{pf.lem.sol.Ialbe.ab.lb=}%
\end{equation}

Write the $\left(  n+m\right)  $-tuple $\mathbf{c}$ in the form $\mathbf{c}%
=\left(  c_{1},c_{2},\ldots,c_{n+m}\right)  $. Thus,%
\[
\left(  c_{1},c_{2},\ldots,c_{n+m}\right)  =\mathbf{c}=\left(  a_{1}%
,a_{2},\ldots,a_{n},b_{1},b_{2},\ldots,b_{m}\right)
\]
(by the definition of $\mathbf{c}$). In other words,%
\begin{equation}
c_{i}=%
\begin{cases}
a_{i}, & \text{if }i\leq n;\\
b_{i-n}, & \text{if }i>n
\end{cases}
\ \ \ \ \ \ \ \ \ \ \text{for every }i\in\left\{  1,2,\ldots,n+m\right\}  .
\label{pf.lem.sol.Ialbe.ab.ci.1}%
\end{equation}
Now, we have%
\begin{equation}
c_{i}=a_{i}\ \ \ \ \ \ \ \ \ \ \text{for every }i\in\left\{  1,2,\ldots
,n\right\}  \label{pf.lem.sol.Ialbe.ab.ci.a}%
\end{equation}
\footnote{\textit{Proof of (\ref{pf.lem.sol.Ialbe.ab.ci.a}):} Let
$i\in\left\{  1,2,\ldots,n\right\}  $. Then, $1\leq i\leq n$. Now,
$n+\underbrace{m}_{\geq0}\geq n$, so that $n\leq n+m$ and thus $i\leq n\leq
n+m$. Combining this with $1\leq i$, we obtain $1\leq i\leq n+m$, so that
$i\in\left\{  1,2,\ldots,n+m\right\}  $. Thus, (\ref{pf.lem.sol.Ialbe.ab.ci.1}%
) yields $c_{i}=%
\begin{cases}
a_{i}, & \text{if }i\leq n;\\
b_{i-n}, & \text{if }i>n
\end{cases}
=a_{i}$ (since $i\leq n$). This proves (\ref{pf.lem.sol.Ialbe.ab.ci.a}).}.
Also,%
\begin{equation}
c_{i}=b_{i-n}\ \ \ \ \ \ \ \ \ \ \text{for every }i\in\left\{  n+1,n+2,\ldots
,n+m\right\}  \label{pf.lem.sol.Ialbe.ab.ci.b}%
\end{equation}
\footnote{\textit{Proof of (\ref{pf.lem.sol.Ialbe.ab.ci.a}):} Let
$i\in\left\{  n+1,n+2,\ldots,n+m\right\}  $. Then, $n+1\leq i\leq n+m$. Now,
$\underbrace{n}_{\geq0}+1\geq1$, so that $1\leq n+1\leq i$. Combining this
with $i\leq n+m$, we obtain $1\leq i\leq n+m$, so that $i\in\left\{
1,2,\ldots,n+m\right\}  $. Thus, (\ref{pf.lem.sol.Ialbe.ab.ci.1}) yields
$c_{i}=%
\begin{cases}
a_{i}, & \text{if }i\leq n;\\
b_{i-n}, & \text{if }i>n
\end{cases}
=b_{i-n}$ (since $i>n$). This proves (\ref{pf.lem.sol.Ialbe.ab.ci.b}).}.

But recall that $\mathbf{c}=\left(  c_{1},c_{2},\ldots,c_{n+m}\right)  $ is an
$\left(  n+m\right)  $-tuple of integers. Thus, Lemma
\ref{lem.sol.Ialbe.iverson} (applied to $n+m$, $\mathbf{c}$ and $c_{i}$
instead of $n$, $\mathbf{a}$ and $a_{i}$) yields%
\begin{align}
\ell\left(  \mathbf{c}\right)   &  =\sum_{i=1}^{n+m}\sum_{j=1}^{n+m}\left[
i<j\right]  \left[  c_{i}>c_{j}\right] \nonumber\\
&  =\sum_{i=1}^{n}\sum_{j=1}^{n+m}\left[  i<j\right]  \left[  c_{i}%
>c_{j}\right]  +\sum_{i=n+1}^{n+m}\sum_{j=1}^{n+m}\left[  i<j\right]  \left[
c_{i}>c_{j}\right]  \label{pf.lem.sol.Ialbe.ab.1}%
\end{align}
(since $0\leq n\leq n+m$). But every $i\in\left\{  1,2,\ldots,n\right\}  $
satisfies%
\begin{equation}
\sum_{j=1}^{n+m}\left[  i<j\right]  \left[  c_{i}>c_{j}\right]  =\sum
_{j=1}^{n}\left[  i<j\right]  \left[  a_{i}>a_{j}\right]  +\sum_{j=1}%
^{m}\left[  a_{i}>b_{j}\right]  \label{pf.lem.sol.Ialbe.ab.3a}%
\end{equation}
\footnote{\textit{Proof of (\ref{pf.lem.sol.Ialbe.ab.3a}):} Let $i\in\left\{
1,2,\ldots,n\right\}  $. Thus, $c_{i}=a_{i}$ (by
(\ref{pf.lem.sol.Ialbe.ab.ci.a})).
\par
For every $j\in\left\{  1,2,\ldots,n\right\}  $, we have%
\begin{equation}
c_{j}=a_{j} \label{pf.lem.sol.Ialbe.ab.3a.pf.1}%
\end{equation}
(by (\ref{pf.lem.sol.Ialbe.ab.ci.a}), applied to $j$ instead of $i$).
\par
For every $j\in\left\{  n+1,n+2,\ldots,n+m\right\}  $, we have%
\begin{equation}
c_{j}=b_{j-n} \label{pf.lem.sol.Ialbe.ab.3a.pf.2}%
\end{equation}
(by (\ref{pf.lem.sol.Ialbe.ab.ci.b}), applied to $j$ instead of $i$).
\par
We have $i\leq n$ (since $i\in\left\{  1,2,\ldots,n\right\}  $). For every
$j\in\left\{  n+1,n+2,\ldots,n+m\right\}  $, we have $j\geq n+1>n\geq i$
(since $i\leq n$). Hence, for every $j\in\left\{  n+1,n+2,\ldots,n+m\right\}
$, we have $i<j$. Thus, for every $j\in\left\{  n+1,n+2,\ldots,n+m\right\}  $,
we have%
\begin{equation}
\left[  i<j\right]  =1. \label{pf.lem.sol.Ialbe.ab.3a.pf.3}%
\end{equation}
\par
Recall that $0\leq n\leq n+m$. Hence,%
\begin{align*}
&  \sum_{j=1}^{n+m}\left[  i<j\right]  \left[  c_{i}>c_{j}\right] \\
&  =\sum_{j=1}^{n}\left[  i<j\right]  \left[  \underbrace{c_{i}}_{=a_{i}%
}>\underbrace{c_{j}}_{\substack{=a_{j}\\\text{(by
(\ref{pf.lem.sol.Ialbe.ab.3a.pf.1}))}}}\right]  +\sum_{j=n+1}^{n+m}%
\underbrace{\left[  i<j\right]  }_{\substack{=1\\\text{(by
(\ref{pf.lem.sol.Ialbe.ab.3a.pf.3}))}}}\left[  \underbrace{c_{i}}_{=a_{i}%
}>\underbrace{c_{j}}_{\substack{=b_{j-n}\\\text{(by
(\ref{pf.lem.sol.Ialbe.ab.3a.pf.2}))}}}\right] \\
&  =\sum_{j=1}^{n}\left[  i<j\right]  \left[  a_{i}>a_{j}\right]
+\sum_{j=n+1}^{n+m}\left[  a_{i}>b_{j-n}\right] \\
&  =\sum_{j=1}^{n}\left[  i<j\right]  \left[  a_{i}>a_{j}\right]  +\sum
_{j=1}^{m}\left[  a_{i}>b_{j}\right] \\
&  \ \ \ \ \ \ \ \ \ \ \left(  \text{here, we have substituted }j\text{ for
}j-n\text{ in the second sum}\right)  .
\end{align*}
This proves (\ref{pf.lem.sol.Ialbe.ab.3a}).}. Also, every $i\in\left\{
n+1,n+2,\ldots,n+m\right\}  $ satisfies%
\begin{equation}
\sum_{j=1}^{n+m}\left[  i<j\right]  \left[  c_{i}>c_{j}\right]  =\sum
_{j=1}^{m}\left[  i-n<j\right]  \left[  b_{i-n}>b_{j}\right]
\label{pf.lem.sol.Ialbe.ab.3b}%
\end{equation}
\footnote{\textit{Proof of (\ref{pf.lem.sol.Ialbe.ab.3b}):} Let $i\in\left\{
n+1,n+2,\ldots,n+m\right\}  $. Thus, $c_{i}=b_{i-n}$ (by
(\ref{pf.lem.sol.Ialbe.ab.ci.b})).
\par
For every $j\in\left\{  n+1,n+2,\ldots,n+m\right\}  $, we have%
\begin{equation}
c_{j}=b_{j-n} \label{pf.lem.sol.Ialbe.ab.3b.pf.2}%
\end{equation}
(by (\ref{pf.lem.sol.Ialbe.ab.ci.b}), applied to $j$ instead of $i$).
\par
We have $i\geq n+1$ (since $i\in\left\{  n+1,n+2,\ldots,n+m\right\}  $) and
thus $i\geq n+1>n$, so that $n<i$. For every $j\in\left\{  1,2,\ldots
,n\right\}  $, we have $j\leq n\leq i$ (since $i\geq n$) and thus $i\geq j$.
Hence, for every $j\in\left\{  1,2,\ldots,n\right\}  $, we don't have $i<j$.
Thus, for every $j\in\left\{  1,2,\ldots,n\right\}  $, we have%
\begin{equation}
\left[  i<j\right]  =0. \label{pf.lem.sol.Ialbe.ab.3b.pf.3}%
\end{equation}
\par
Recall that $0\leq n\leq n+m$. Hence,%
\begin{align*}
&  \sum_{j=1}^{n+m}\left[  i<j\right]  \left[  c_{i}>c_{j}\right] \\
&  =\sum_{j=1}^{n}\underbrace{\left[  i<j\right]  }_{\substack{=0\\\text{(by
(\ref{pf.lem.sol.Ialbe.ab.3b.pf.3}))}}}\left[  c_{i}>c_{j}\right]
+\sum_{j=n+1}^{n+m}\underbrace{\left[  i<j\right]  }_{\substack{=\left[
i-n<j-n\right]  \\\text{(since }i<j\text{ is equivalent to}\\i-n<j-n\text{)}%
}}\left[  \underbrace{c_{i}}_{=b_{i-n}}>\underbrace{c_{j}}_{\substack{=b_{j-n}%
\\\text{(by (\ref{pf.lem.sol.Ialbe.ab.3b.pf.2}))}}}\right] \\
&  =\underbrace{\sum_{j=1}^{n}0\left[  c_{i}>c_{j}\right]  }_{=0}+\sum
_{j=n+1}^{n+m}\left[  i-n<j-n\right]  \left[  b_{i-n}>b_{j-n}\right] \\
&  =\sum_{j=n+1}^{n+m}\left[  i-n<j-n\right]  \left[  b_{i-n}>b_{j-n}\right]
\\
&  =\sum_{j=1}^{m}\left[  i-n<j\right]  \left[  b_{i-n}>b_{j}\right] \\
&  \ \ \ \ \ \ \ \ \ \ \left(  \text{here, we have substituted }j\text{ for
}j-n\text{ in the sum}\right)  .
\end{align*}
This proves (\ref{pf.lem.sol.Ialbe.ab.3b}).}.

Now, (\ref{pf.lem.sol.Ialbe.ab.1}) becomes%
\begin{align*}
\ell\left(  \mathbf{c}\right)   &  =\sum_{i=1}^{n}\underbrace{\sum_{j=1}%
^{n+m}\left[  i<j\right]  \left[  c_{i}>c_{j}\right]  }_{\substack{=\sum
_{j=1}^{n}\left[  i<j\right]  \left[  a_{i}>a_{j}\right]  +\sum_{j=1}%
^{m}\left[  a_{i}>b_{j}\right]  \\\text{(by (\ref{pf.lem.sol.Ialbe.ab.3a}))}%
}}\\
&  \ \ \ \ \ \ \ \ \ \ +\sum_{i=n+1}^{n+m}\underbrace{\sum_{j=1}^{n+m}\left[
i<j\right]  \left[  c_{i}>c_{j}\right]  }_{\substack{=\sum_{j=1}^{m}\left[
i-n<j\right]  \left[  b_{i-n}>b_{j}\right]  \\\text{(by
(\ref{pf.lem.sol.Ialbe.ab.3b}))}}}\\
&  =\underbrace{\sum_{i=1}^{n}\left(  \sum_{j=1}^{n}\left[  i<j\right]
\left[  a_{i}>a_{j}\right]  +\sum_{j=1}^{m}\left[  a_{i}>b_{j}\right]
\right)  }_{=\sum_{i=1}^{n}\sum_{j=1}^{n}\left[  i<j\right]  \left[
a_{i}>a_{j}\right]  +\sum_{i=1}^{n}\sum_{j=1}^{m}\left[  a_{i}>b_{j}\right]
}\\
&  \ \ \ \ \ \ \ \ \ \ +\underbrace{\sum_{i=n+1}^{n+m}\sum_{j=1}^{m}\left[
i-n<j\right]  \left[  b_{i-n}>b_{j}\right]  }_{\substack{=\sum_{i=1}^{m}%
\sum_{j=1}^{m}\left[  i<j\right]  \left[  b_{i}>b_{j}\right]  \\\text{(here,
we have substituted }i\text{ for }i-n\\\text{in the outer sum)}}}\\
&  =\underbrace{\sum_{i=1}^{n}\sum_{j=1}^{n}\left[  i<j\right]  \left[
a_{i}>a_{j}\right]  }_{\substack{=\ell\left(  \mathbf{a}\right)  \\\text{(by
(\ref{pf.lem.sol.Ialbe.ab.la=}))}}}+\underbrace{\sum_{i=1}^{n}}%
_{\substack{=\sum_{i\in\left\{  1,2,\ldots,n\right\}  }=\sum_{i\in\left[
n\right]  }\\\text{(since }\left\{  1,2,\ldots,n\right\}  =\left[  n\right]
\\\text{(because }\left[  n\right]  =\left\{  1,2,\ldots,n\right\}
\\\text{(by the definition of }\left[  n\right]  \text{)))}}}\underbrace{\sum
_{j=1}^{m}}_{\substack{=\sum_{j\in\left\{  1,2,\ldots,m\right\}  }=\sum
_{j\in\left[  m\right]  }\\\text{(since }\left\{  1,2,\ldots,m\right\}
=\left[  m\right]  \\\text{(because }\left[  m\right]  =\left\{
1,2,\ldots,m\right\}  \\\text{(by the definition of }\left[  m\right]
\text{)))}}}\left[  a_{i}>b_{j}\right] \\
&  \ \ \ \ \ \ \ \ \ \ +\underbrace{\sum_{i=1}^{m}\sum_{j=1}^{m}\left[
i<j\right]  \left[  b_{i}>b_{j}\right]  }_{\substack{=\ell\left(
\mathbf{b}\right)  \\\text{(by (\ref{pf.lem.sol.Ialbe.ab.lb=}))}}}\\
&  =\ell\left(  \mathbf{a}\right)  +\underbrace{\sum_{i\in\left[  n\right]
}\sum_{j\in\left[  m\right]  }}_{=\sum_{\left(  i,j\right)  \in\left[
n\right]  \times\left[  m\right]  }}\left[  a_{i}>b_{j}\right]  +\ell\left(
\mathbf{b}\right) \\
&  =\ell\left(  \mathbf{a}\right)  +\sum_{\left(  i,j\right)  \in\left[
n\right]  \times\left[  m\right]  }\left[  a_{i}>b_{j}\right]  +\ell\left(
\mathbf{b}\right)  =\ell\left(  \mathbf{a}\right)  +\ell\left(  \mathbf{b}%
\right)  +\sum_{\left(  i,j\right)  \in\left[  n\right]  \times\left[
m\right]  }\left[  a_{i}>b_{j}\right]  .
\end{align*}
This proves Lemma \ref{lem.sol.Ialbe.ab}.
\end{proof}
\end{verlong}

\begin{lemma}
\label{lem.sol.Ialbe.inclist}Let $S$ be a finite set of integers. Let $\left(
c_{1},c_{2},\ldots,c_{s}\right)  $ be a list of all elements of $S$ (with no repetitions).

\textbf{(a)} Then, the map $\left[  s\right]  \rightarrow S,\ h\mapsto c_{h}$
is well-defined and a bijection.

\textbf{(b)} Let $\pi\in S_{s}$. Then, $\left(  c_{\pi\left(  1\right)
},c_{\pi\left(  2\right)  },\ldots,c_{\pi\left(  s\right)  }\right)  $ is a
list of all elements of $S$ (with no repetitions).
\end{lemma}

\begin{vershort}
\begin{proof}
[Proof of Lemma \ref{lem.sol.Ialbe.inclist}.]Lemma \ref{lem.sol.Ialbe.inclist}
should really be intuitively obvious; let me merely sketch how to formalize
the intuition.

We have the following basic fact:

\begin{statement}
\textit{Statement 1:} Let $X$ and $Y$ be two sets. Let $\phi:X\rightarrow Y$
be a map. Then, $\phi$ is a bijection if and only if each $i\in Y$ has exactly
one preimage under $\phi$.
\end{statement}

Statement 1 is well-known and easy to prove; it is a rather useful (necessary
and sufficient) criterion for proving the bijectivity of maps (particularly
since it does not require proving injectivity and surjectivity separately).

We assumed that $\left(  c_{1},c_{2},\ldots,c_{s}\right)  $ is a list of all
elements of $S$ (with no repetitions). This means that the following two
statements are valid:

\begin{statement}
\textit{Statement 2:} The list $\left(  c_{1},c_{2},\ldots,c_{s}\right)  $ is
a list of elements of $S$. (In other words, we have $c_{h}\in S$ for every
$h\in\left[  s\right]  $.)
\end{statement}

\begin{statement}
\textit{Statement 3:} Each element of $S$ appears exactly once in the list
$\left(  c_{1},c_{2},\ldots,c_{s}\right)  $. In other words, for each $i\in
S$, there exists exactly one $h\in\left[  s\right]  $ satisfying $i=c_{h}$.
\end{statement}

Statement 2 shows that $c_{h}\in S$ for every $h\in\left[  s\right]  $. Thus,
the map $\left[  s\right]  \rightarrow S,\ h\mapsto c_{h}$ is well-defined.
Denote this map by $\alpha$. Thus, $\alpha\left(  h\right)  =c_{h}$ for every
$h\in\left[  s\right]  $.

\textbf{(a)} Statement 3 shows that, for each $i\in S$, there exists exactly
one $h\in\left[  s\right]  $ satisfying $i=c_{h}$. Since $\alpha\left(
h\right)  =c_{h}$ for every $h\in\left[  s\right]  $, we can rewrite this as
follows: For each $i\in S$, there exists exactly one $h\in\left[  s\right]  $
satisfying $i=\alpha\left(  h\right)  $. In other words, each $i\in S$ has
exactly one preimage under $\alpha$. According to Statement 1 (applied to
$X=\left[  s\right]  $, $Y=S$ and $\phi=\alpha$), this holds if and only if
$\alpha$ is a bijection. Thus, $\alpha$ is a bijection. In other words, the
map $\left[  s\right]  \rightarrow S,\ h\mapsto c_{h}$ is a bijection (because
$\alpha$ is precisely this map). Thus, Lemma \ref{lem.sol.Ialbe.inclist}
\textbf{(a)} is proven.

\textbf{(b)} We must prove that $\left(  c_{\pi\left(  1\right)  }%
,c_{\pi\left(  2\right)  },\ldots,c_{\pi\left(  s\right)  }\right)  $ is a
list of all elements of $S$ (with no repetitions). This means proving the
following two statements:

\begin{statement}
\textit{Statement 4:} The list $\left(  c_{\pi\left(  1\right)  }%
,c_{\pi\left(  2\right)  },\ldots,c_{\pi\left(  s\right)  }\right)  $ is a
list of elements of $S$. (In other words, we have $c_{\pi\left(  h\right)
}\in S$ for every $h\in\left[  s\right]  $.)
\end{statement}

\begin{statement}
\textit{Statement 5:} Each element of $S$ appears exactly once in the list
\newline$\left(  c_{\pi\left(  1\right)  },c_{\pi\left(  2\right)  }%
,\ldots,c_{\pi\left(  s\right)  }\right)  $. In other words, for each $i\in
S$, there exists exactly one $h\in\left[  s\right]  $ satisfying
$i=c_{\pi\left(  h\right)  }$.
\end{statement}

[\textit{Proof of Statement 4:} Statement 2 shows that $c_{h}\in S$ for every
$h\in\left[  s\right]  $. Applying this to $\pi\left(  h\right)  $ instead of
$h$, we conclude that $c_{\pi\left(  h\right)  }\in S$ for every $h\in\left[
s\right]  $. This proves Statement 4.]

[\textit{Proof of Statement 5:} The map $\alpha\circ\pi$ is a bijection (since
it is the composition of the two bijections $\alpha$ and $\pi$). According to
Statement 1 (applied to $X=\left[  s\right]  $, $Y=S$ and $\phi=\alpha\circ
\pi$), this means that each $i\in S$ has exactly one preimage under
$\alpha\circ\pi$. In other words, for each $i\in S$, there exists exactly one
$h\in\left[  s\right]  $ satisfying $i=\left(  \alpha\circ\pi\right)  \left(
h\right)  $. Since every $h\in\left[  s\right]  $ satisfies%
\[
\left(  \alpha\circ\pi\right)  \left(  h\right)  =\alpha\left(  \pi\left(
h\right)  \right)  =c_{\pi\left(  h\right)  }\ \ \ \ \ \ \ \ \ \ \left(
\text{by the definition of }\alpha\right)  ,
\]
we can rewrite this as follows: For each $i\in S$, there exists exactly one
$h\in\left[  s\right]  $ satisfying $i=c_{\pi\left(  h\right)  }$. This proves
Statement 5.]

Now, Statements 4 and 5 are proven; thus, $\left(  c_{\pi\left(  1\right)
},c_{\pi\left(  2\right)  },\ldots,c_{\pi\left(  s\right)  }\right)  $ is a
list of all elements of $S$ (with no repetitions). This proves Lemma
\ref{lem.sol.Ialbe.inclist} \textbf{(b)}.
\end{proof}
\end{vershort}

\begin{verlong}
\begin{proof}
[Proof of Lemma \ref{lem.sol.Ialbe.inclist}.]We have $\left[  s\right]
=\left\{  1,2,\ldots,s\right\}  $ (by the definition of $\left[  s\right]  $).

The list $\left(  c_{1},c_{2},\ldots,c_{s}\right)  $ be a list of all elements
of $S$ (with no repetitions). In other words, the list $\left(  c_{1}%
,c_{2},\ldots,c_{s}\right)  $ is a list of elements of $S$, and contains each
element of $S$ exactly once.

The list $\left(  c_{1},c_{2},\ldots,c_{s}\right)  $ contains each element of
$S$ exactly once. In other words, for every $i\in S$, the list $\left(
c_{1},c_{2},\ldots,c_{s}\right)  $ contains the element $i$ exactly once. In
other words,
\begin{equation}
\text{for every }i\in S\text{, there exists a unique }p\in\left\{
1,2,\ldots,s\right\}  \text{ such that }i=c_{p}.
\label{pf.lem.sol.Ialbe.inclist.exonce}%
\end{equation}

For every $h\in\left[  s\right]  $, the element $c_{h}$ of $S$ is
well-defined\footnote{\textit{Proof.} Let $h\in\left[  s\right]  $. Then,
$h\in\left[  s\right]  =\left\{  1,2,\ldots,s\right\}  $. Thus, the element
$c_{h}$ of $S$ is well-defined (since $\left(  c_{1},c_{2},\ldots
,c_{s}\right)  $ is a list of elements of $S$). Qed.}. Thus, the map $\left[
s\right]  \rightarrow S,\ h\mapsto c_{h}$ is well-defined. Let us denote this
map by $\phi$. Thus, $\phi$ is the map $\left[  s\right]  \rightarrow
S,\ h\mapsto c_{h}$.

The map $\phi$ is injective\footnote{\textit{Proof.} Let $g$ and $h$ be two
elements of $\left[  s\right]  $ such that $\phi\left(  g\right)  =\phi\left(
h\right)  $. We shall show that $g=h$.
\par
We have $g\in\left[  s\right]  =\left\{  1,2,\ldots,s\right\}  $ and
$h\in\left[  s\right]  =\left\{  1,2,\ldots,s\right\}  $.
\par
The definition of $\phi$ yields $\phi\left(  g\right)  =c_{g}$. Hence,
$c_{g}=\phi\left(  g\right)  =\phi\left(  h\right)  =c_{h}$ (by the definition
of $\phi$). Hence, $h$ is an element of $\left\{  1,2,\ldots,s\right\}  $
satisfying $c_{g}=c_{h}$. In other words, $h$ is an element $p\in\left\{
1,2,\ldots,s\right\}  $ satisfying $c_{g}=c_{p}$.
\par
Also, $g$ is an element of $\left\{  1,2,\ldots,s\right\}  $ satisfying
$c_{g}=c_{g}$. In other words, $g$ is an element $p\in\left\{  1,2,\ldots
,s\right\}  $ satisfying $c_{g}=c_{p}$.
\par
But $c_{g}\in S$ (since $\left(  c_{1},c_{2},\ldots,c_{s}\right)  $ is a list
of elements of $S$). Hence, (\ref{pf.lem.sol.Ialbe.inclist.exonce}) (applied
to $i=c_{g}$) shows that there exists a unique $p\in\left\{  1,2,\ldots
,s\right\}  $ such that $c_{g}=c_{p}$. In particular, there exists \textbf{at
most one} such $p$. In other words, if $u$ and $v$ are two elements
$p\in\left\{  1,2,\ldots,s\right\}  $ satisfying $c_{g}=c_{p}$, then $u=v$.
Applying this to $u=g$ and $v=h$, we obtain $g=h$ (since $g$ and $h$ are two
elements $p\in\left\{  1,2,\ldots,s\right\}  $ satisfying $c_{g}=c_{p}$).
\par
Now, forget that we fixed $g$ and $h$. We thus have proven that if $g$ and $h$
are two elements of $\left[  s\right]  $ such that $\phi\left(  g\right)
=\phi\left(  h\right)  $, then $g=h$. In other words, the map $\phi$ is
injective. Qed.} and surjective\footnote{\textit{Proof.} Let $i\in S$. Then,
there exists a unique $p\in\left\{  1,2,\ldots,s\right\}  $ such that
$i=c_{p}$ (by (\ref{pf.lem.sol.Ialbe.inclist.exonce})). Consider this $p$. We
have $p\in\left\{  1,2,\ldots,s\right\}  =\left[  s\right]  $ (since $\left[
s\right]  =\left\{  1,2,\ldots,s\right\}  $). Thus, $\phi\left(  p\right)  $
is well-defined. Now, the definition of $\phi$ yields $\phi\left(  p\right)
=c_{p}=i$. Thus, $i=\phi\left(  \underbrace{p}_{\in\left[  s\right]  }\right)
\in\phi\left(  \left[  s\right]  \right)  $.
\par
Now, forget that we fixed $i$. We thus have shown that every $i\in S$
satisfies $i\in\phi\left(  \left[  s\right]  \right)  $. In other words,
$S\subseteq\phi\left(  \left[  s\right]  \right)  $. In other words, the map
$\phi$ is surjective. Qed.}. Hence, this map $\phi$ is bijective. In other
words, $\phi$ is a bijection. In other words, the map $\left[  s\right]
\rightarrow S,\ h\mapsto c_{h}$ is a bijection (since $\phi$ is the map
$\left[  s\right]  \rightarrow S,\ h\mapsto c_{h}$). This completes the proof
of Lemma \ref{lem.sol.Ialbe.inclist} \textbf{(a)}.

\textbf{(b)} We know that $S_{s}$ is the set of all permutations of $\left\{
1,2,\ldots,s\right\}  $ (by the definition of $S_{s}$). In other words,
$S_{s}$ is the set of all permutations of $\left[  s\right]  $ (since
$\left\{  1,2,\ldots,s\right\}  =\left[  s\right]  $). Now, $\pi$ is an
element of $S_{s}$. In other words, $\pi$ is a permutation of $\left[
s\right]  $ (since $S_{s}$ is the set of all permutations of $\left[
s\right]  $). In other words, $\pi$ is a bijection $\left[  s\right]
\rightarrow\left[  s\right]  $. So the map $\pi:\left[  s\right]
\rightarrow\left[  s\right]  $ is bijective. Thus, this map $\pi$ is
surjective and injective.

The map $\pi$ is a map $\left[  s\right]  \rightarrow\left[  s\right]  $. In
other words, $\pi$ is a map $\left\{  1,2,\ldots,s\right\}  \rightarrow
\left\{  1,2,\ldots,s\right\}  $ (since $\left[  s\right]  =\left\{
1,2,\ldots,s\right\}  $). Hence, for every $i\in\left\{  1,2,\ldots,s\right\}
$, we have $\pi\left(  i\right)  \in\left\{  1,2,\ldots,s\right\}  $. Thus,
for every $i\in\left\{  1,2,\ldots,s\right\}  $, the element $c_{\pi\left(
i\right)  }$ is a well-defined element of $S$. Thus, $\left(  c_{\pi\left(
1\right)  },c_{\pi\left(  2\right)  },\ldots,c_{\pi\left(  s\right)  }\right)
$ is a list of elements of $S$.

\begin{noncompile}
The map $\phi\circ\pi:\left[  s\right]  \rightarrow S$ is bijective (since it
is the composition of the two bijective maps $\phi$ and $\pi$), and thus
injective and surjective.
\end{noncompile}

Let $i\in S$. We shall show that there exists a unique $p\in\left[  s\right]
$ such that $i=c_{\pi\left(  p\right)  }$.

If $g$ and $h$ are two elements $p\in\left[  s\right]  $ such that
$i=c_{\pi\left(  p\right)  }$, then $g=h$\ \ \ \ \footnote{\textit{Proof.} Let
$g$ and $h$ be two elements $p\in\left[  s\right]  $ such that $i=c_{\pi
\left(  p\right)  }$. We must prove that $g=h$.
\par
We know that $g$ is an element $p\in\left[  s\right]  $ such that
$i=c_{\pi\left(  p\right)  }$. In other words, $g$ is an element of $\left[
s\right]  $ and satisfies $i=c_{\pi\left(  g\right)  }$. We have $\pi\left(
g\right)  \in\left[  s\right]  $ (since $\pi$ is a map $\left[  s\right]
\rightarrow\left[  s\right]  $) and thus $\phi\left(  \pi\left(  g\right)
\right)  =c_{\pi\left(  g\right)  }$ (by the definition of $\phi$). Comparing
this with $i=c_{\pi\left(  g\right)  }$, we obtain $\phi\left(  \pi\left(
g\right)  \right)  =i$.
\par
We know that $h$ is an element $p\in\left[  s\right]  $ such that
$i=c_{\pi\left(  p\right)  }$. In other words, $h$ is an element of $\left[
s\right]  $ and satisfies $i=c_{\pi\left(  h\right)  }$. We have $\pi\left(
h\right)  \in\left[  s\right]  $ (since $\pi$ is a map $\left[  s\right]
\rightarrow\left[  s\right]  $) and thus $\phi\left(  \pi\left(  h\right)
\right)  =c_{\pi\left(  h\right)  }$ (by the definition of $\phi$). Comparing
this with $i=c_{\pi\left(  h\right)  }$, we obtain $\phi\left(  \pi\left(
h\right)  \right)  =i$.
\par
Now, $\phi\left(  \pi\left(  g\right)  \right)  =i=\phi\left(  \pi\left(
h\right)  \right)  $. This yields $\pi\left(  g\right)  =\pi\left(  h\right)
$ (since the map $\phi$ is injective). Therefore, $g=h$ (since the map $\pi$
is injective). Qed.}. In other words,
\begin{equation}
\text{there exists \textbf{at most one} }p\in\left[  s\right]  \text{ such
that }i=c_{\pi\left(  p\right)  }. \label{pf.lem.sol.Ialbe.inclist.b.1}%
\end{equation}

On the other hand, $i\in S=\phi\left(  \left[  s\right]  \right)  $ (since the
map $\phi$ is surjective). In other words, there exists some $h\in\left[
s\right]  $ such that $i=\phi\left(  h\right)  $. Consider this $h$. We have
$h\in\left[  s\right]  =\pi\left(  \left[  s\right]  \right)  $ (since the map
$\pi$ is surjective). In other words, there exists some $g\in\left[  s\right]
$ such that $h=\pi\left(  g\right)  $. Consider this $g$.

Now, the definition of $\phi$ yields $\phi\left(  h\right)  =c_{h}$, so that
$c_{h}=\phi\left(  h\right)  =i$. Thus, $i=c_{h}=c_{\pi\left(  g\right)  }$
(since $h=\pi\left(  g\right)  $). Thus, $g$ is an element of $\left[
s\right]  $ and satisfies $i=c_{\pi\left(  g\right)  }$. In other words, $g$
is a $p\in\left[  s\right]  $ such that $i=c_{\pi\left(  p\right)  }$. Hence,
there exists \textbf{at least one} $p\in\left[  s\right]  $ such that
$i=c_{\pi\left(  p\right)  }$ (namely, $p=g$). Combining this with
(\ref{pf.lem.sol.Ialbe.inclist.b.1}), we conclude that there exists a unique
$p\in\left[  s\right]  $ such that $i=c_{\pi\left(  p\right)  }$. In other
words, there exists a unique $p\in\left\{  1,2,\ldots,s\right\}  $ such that
$i=c_{\pi\left(  p\right)  }$ (since $\left[  s\right]  =\left\{
1,2,\ldots,s\right\}  $).

Now, forget that we fixed $i$. We thus have shown that%
\begin{equation}
\text{for every }i\in S\text{, there exists a unique }p\in\left\{
1,2,\ldots,s\right\}  \text{ such that }i=c_{\pi\left(  p\right)  }.
\label{pf.lem.sol.Ialbe.inclist.b.2}%
\end{equation}
In other words, the list $\left(  c_{\pi\left(  1\right)  },c_{\pi\left(
2\right)  },\ldots,c_{\pi\left(  s\right)  }\right)  $ contains each element
of $S$ exactly once. Since $\left(  c_{\pi\left(  1\right)  },c_{\pi\left(
2\right)  },\ldots,c_{\pi\left(  s\right)  }\right)  $ is a list of elements
of $S$, this yields that\newline$\left(  c_{\pi\left(  1\right)  }%
,c_{\pi\left(  2\right)  },\ldots,c_{\pi\left(  s\right)  }\right)  $ is a
list of all elements of $S$ (with no repetitions). This proves Lemma
\ref{lem.sol.Ialbe.inclist} \textbf{(b)}.
\end{proof}
\end{verlong}

\begin{lemma}
\label{lem.sol.Ialbe.II}Let $I$ be a finite set of integers. Let $k=\left\vert
I\right\vert $. Then,%
\[
\sum_{x\in I}\sum_{y\in I}\left[  x>y\right]  =0+1+\cdots+\left(  k-1\right)
.
\]

\end{lemma}

\begin{vershort}
\begin{proof}
[Proof of Lemma \ref{lem.sol.Ialbe.II}.]Let $\left(  i_{1},i_{2},\ldots
,i_{k}\right)  $ be the list of all elements of $I$ in increasing order (with
no repetitions). (Such a list exists, since $\left\vert I\right\vert =k$.)
Then, Lemma \ref{lem.sol.Ialbe.inclist} \textbf{(a)} (applied to $I$, $k$ and
$\left(  i_{1},i_{2},\ldots,i_{k}\right)  $ instead of $S$, $s$ and $\left(
c_{1},c_{2},\ldots,c_{s}\right)  $) shows that the map $\left[  k\right]
\rightarrow I,\ h\mapsto i_{h}$ is well-defined and a bijection.

We have $i_{1}<i_{2}<\cdots<i_{k}$ (because of how $\left(  i_{1},i_{2}%
,\ldots,i_{k}\right)  $ was defined). Hence, if $u$ and $v$ are two elements
of $\left[  k\right]  $, then $i_{u}>i_{v}$ holds if and only if $u>v$. In
other words, if $u$ and $v$ are two elements of $\left[  k\right]  $, then%
\begin{equation}
\left[  i_{u}>i_{v}\right]  =\left[  u>v\right]  .
\label{pf.lem.sol.Ialbe.II.short.3}%
\end{equation}

Now, every $x\in I$ satisfies $\sum_{y\in I}\left[  x>y\right]  =\sum
_{v\in\left[  k\right]  }\left[  x>i_{v}\right]  $ (here, we have substituted
$i_{v}$ for $y$ in the sum, since the map $\left[  k\right]  \rightarrow
I,\ h\mapsto i_{h}$ is a bijection). Thus,%
\begin{align}
\sum_{x\in I}\underbrace{\sum_{y\in I}\left[  x>y\right]  }_{=\sum
_{v\in\left[  k\right]  }\left[  x>i_{v}\right]  }  &  =\sum_{x\in I}%
\sum_{v\in\left[  k\right]  }\left[  x>i_{v}\right] \nonumber\\
&  =\sum_{u\in\left[  k\right]  }\sum_{v\in\left[  k\right]  }%
\underbrace{\left[  i_{u}>i_{v}\right]  }_{\substack{=\left[  u>v\right]
\\\text{(by (\ref{pf.lem.sol.Ialbe.II.short.3}))}}}\ \ \ \ \ \ \ \ \ \ \left(
\begin{array}
[c]{c}%
\text{here, we have substituted }i_{u}\\
\text{for }x\text{ in the sum, since the}\\
\text{map }\left[  k\right]  \rightarrow I,\ h\mapsto i_{h}\text{ is a
bijection}%
\end{array}
\right) \nonumber\\
&  =\sum_{u\in\left[  k\right]  }\sum_{v\in\left[  k\right]  }\left[
u>v\right]  . \label{pf.lem.sol.Ialbe.II.short.4}%
\end{align}
But every $u\in\left[  k\right]  $ satisfies%
\begin{equation}
\sum_{v\in\left[  k\right]  }\left[  u>v\right]  =u-1
\label{pf.lem.sol.Ialbe.II.short.5}%
\end{equation}
\footnote{\textit{Proof of (\ref{pf.lem.sol.Ialbe.II.short.5}):} Let
$u\in\left[  k\right]  $. Then,%
\begin{align*}
\underbrace{\sum_{v\in\left[  k\right]  }}_{=\sum_{v=1}^{k}}\left[
u>v\right]   &  =\sum_{v=1}^{k}\left[  u>v\right]  =\sum_{v=1}^{u-1}%
\underbrace{\left[  u>v\right]  }_{\substack{=1\\\text{(since }%
u>v\\\text{(since }v\leq u-1<u\text{))}}}+\sum_{v=u}^{k}\underbrace{\left[
u>v\right]  }_{\substack{=0\\\text{(since we don't have }u>v\\\text{(since
}v\geq u\text{))}}}\ \ \ \ \ \ \ \ \ \ \left(  \text{since }1\leq u\leq
k\right) \\
&  =\sum_{v=1}^{u-1}1+\underbrace{\sum_{v=u}^{k}0}_{=0}=\sum_{v=1}%
^{u-1}1=\left(  u-1\right)  1=u-1.
\end{align*}
This proves (\ref{pf.lem.sol.Ialbe.II.short.5}).}.

Now, (\ref{pf.lem.sol.Ialbe.II.short.4}) becomes%
\begin{align*}
\sum_{x\in I}\sum_{y\in I}\left[  x>y\right]   &  =\underbrace{\sum
_{u\in\left[  k\right]  }}_{=\sum_{u=1}^{k}}\underbrace{\sum_{v\in\left[
k\right]  }\left[  u>v\right]  }_{\substack{=u-1\\\text{(by
(\ref{pf.lem.sol.Ialbe.II.short.5}))}}}=\sum_{u=1}^{k}\left(  u-1\right) \\
&  =\sum_{u=0}^{k-1}u\ \ \ \ \ \ \ \ \ \ \left(  \text{here, we have
substituted }u\text{ for }u-1\text{ in the sum}\right) \\
&  =0+1+\cdots+\left(  k-1\right)  .
\end{align*}
This proves Lemma \ref{lem.sol.Ialbe.II}.
\end{proof}
\end{vershort}

\begin{verlong}
\begin{proof}
[Proof of Lemma \ref{lem.sol.Ialbe.II}.]We have $\left[  k\right]  =\left\{
1,2,\ldots,k\right\}  $ (by the definition of $\left[  k\right]  $).

Let $\left(  i_{1},i_{2},\ldots,i_{k}\right)  $ be the list of all elements of
$I$ in increasing order (with no repetitions). (Such a list exists, since
$\left\vert I\right\vert =k$.) Then, Lemma \ref{lem.sol.Ialbe.inclist}
\textbf{(a)} (applied to $I$, $k$ and $\left(  i_{1},i_{2},\ldots
,i_{k}\right)  $ instead of $S$, $s$ and $\left(  c_{1},c_{2},\ldots
,c_{s}\right)  $) shows that the map $\left[  k\right]  \rightarrow
I,\ h\mapsto i_{h}$ is well-defined and a bijection.

We know that $\left(  i_{1},i_{2},\ldots,i_{k}\right)  $ is the list of all
elements of $I$ in increasing order (with no repetitions). Thus, $\left(
i_{1},i_{2},\ldots,i_{k}\right)  $ is a strictly increasing list. In other
words, we have $i_{1}<i_{2}<\cdots<i_{k}$. In other words, if $u$ and $v$ are
two elements of $\left[  k\right]  $ such that $u<v$, then%
\begin{equation}
i_{u}<i_{v}. \label{pf.lem.sol.Ialbe.II.1}%
\end{equation}
Moreover, if $u$ and $v$ are two elements of $\left[  k\right]  $ such that
$u\leq v$, then%
\begin{equation}
i_{u}\leq i_{v} \label{pf.lem.sol.Ialbe.II.2}%
\end{equation}
(since $i_{1}<i_{2}<\cdots<i_{k}$).

Now, every two elements $u$ and $v$ of $\left[  k\right]  $ satisfy
\begin{equation}
\left[  i_{u}>i_{v}\right]  =\left[  u>v\right]  \label{pf.lem.sol.Ialbe.II.3}%
\end{equation}
\footnote{\textit{Proof of (\ref{pf.lem.sol.Ialbe.II.3}):} Let $u$ and $v$ be
two elements of $\left[  k\right]  $. We are in one of the following two
cases:
\par
\textit{Case 1:} We have $u>v$.
\par
\textit{Case 2:} We have $u\leq v$.
\par
Let us first consider Case 1. In this case, we have $u>v$. Hence, $v<u$. Thus,
(\ref{pf.lem.sol.Ialbe.II.1}) (applied to $v$ and $u$ instead of $u$ and $v$)
yields $i_{v}<i_{u}$. In other words, $i_{u}>i_{v}$. Hence, $\left[
i_{u}>i_{v}\right]  =1$. But $u>v$, and thus $\left[  u>v\right]  =1$.
Comparing this with $\left[  i_{u}>i_{v}\right]  =1$, we obtain $\left[
i_{u}>i_{v}\right]  =\left[  u>v\right]  $. Thus, (\ref{pf.lem.sol.Ialbe.II.3}%
) is proven in Case 1.
\par
Let us now consider Case 2. In this case, we have $u\leq v$. Hence, $i_{u}\leq
i_{v}$ (by (\ref{pf.lem.sol.Ialbe.II.2})). In other words, we don't have
$i_{u}>i_{v}$. Hence, $\left[  i_{u}>i_{v}\right]  =0$. But $u\leq v$. In
other words, we don't have $u>v$. Hence, $\left[  u>v\right]  =0$. Comparing
this with $\left[  i_{u}>i_{v}\right]  =0$, we obtain $\left[  i_{u}%
>i_{v}\right]  =\left[  u>v\right]  $. Thus, (\ref{pf.lem.sol.Ialbe.II.3}) is
proven in Case 2.
\par
We thus have proven (\ref{pf.lem.sol.Ialbe.II.3}) in each of the two Cases 1
and 2. Since these two Cases cover all possibilities, this shows that
(\ref{pf.lem.sol.Ialbe.II.3}) always holds. Qed.}.

Now, every $x\in I$ satisfies%
\begin{align*}
\sum_{y\in I}\left[  x>y\right]   &  =\sum_{h\in\left[  k\right]  }\left[
x>i_{h}\right]  \ \ \ \ \ \ \ \ \ \ \left(
\begin{array}
[c]{c}%
\text{here, we have substituted }i_{h}\text{ for }y\text{ in the sum,}\\
\text{since the map }\left[  k\right]  \rightarrow I,\ h\mapsto i_{h}\text{ is
a bijection}%
\end{array}
\right) \\
&  =\sum_{v\in\left[  k\right]  }\left[  x>i_{v}\right]
\ \ \ \ \ \ \ \ \ \ \left(  \text{here, we have renamed the summation index
}h\text{ as }v\right)  .
\end{align*}
Thus,%
\begin{align}
&  \sum_{x\in I}\underbrace{\sum_{y\in I}\left[  x>y\right]  }_{=\sum
_{v\in\left[  k\right]  }\left[  x>i_{v}\right]  }\nonumber\\
&  =\sum_{x\in I}\sum_{v\in\left[  k\right]  }\left[  x>i_{v}\right]
\nonumber\\
&  =\sum_{h\in\left[  k\right]  }\sum_{v\in\left[  k\right]  }\left[
i_{h}>i_{v}\right]  \ \ \ \ \ \ \ \ \ \ \left(
\begin{array}
[c]{c}%
\text{here, we have substituted }i_{h}\text{ for }x\text{ in the outer sum,}\\
\text{since the map }\left[  k\right]  \rightarrow I,\ h\mapsto i_{h}\text{ is
a bijection}%
\end{array}
\right) \nonumber\\
&  =\sum_{u\in\left[  k\right]  }\sum_{v\in\left[  k\right]  }%
\underbrace{\left[  i_{u}>i_{v}\right]  }_{\substack{=\left[  u>v\right]
\\\text{(by (\ref{pf.lem.sol.Ialbe.II.3}))}}}\ \ \ \ \ \ \ \ \ \ \left(
\begin{array}
[c]{c}%
\text{here, we have renamed the summation}\\
\text{index }h\text{ as }u\text{ in the outer sum}%
\end{array}
\right) \nonumber\\
&  =\sum_{u\in\left[  k\right]  }\sum_{v\in\left[  k\right]  }\left[
u>v\right]  . \label{pf.lem.sol.Ialbe.II.4}%
\end{align}
But every $u\in\left[  k\right]  $ satisfies%
\begin{equation}
\sum_{v\in\left[  k\right]  }\left[  u>v\right]  =u-1
\label{pf.lem.sol.Ialbe.II.5}%
\end{equation}
\footnote{\textit{Proof of (\ref{pf.lem.sol.Ialbe.II.5}):} Let $u\in\left[
k\right]  $. Then, $u\in\left[  k\right]  =\left\{  1,2,\ldots,k\right\}  $
(by the definition of $\left[  k\right]  $). Hence, $1\leq u\leq k$. Now,%
\begin{align*}
\underbrace{\sum_{v\in\left[  k\right]  }}_{\substack{=\sum_{v\in\left\{
1,2,\ldots,k\right\}  }\\\text{(since }\left[  k\right]  =\left\{
1,2,\ldots,k\right\}  \text{)}}}\left[  u>v\right]   &  =\underbrace{\sum
_{v\in\left\{  1,2,\ldots,k\right\}  }}_{=\sum_{v=1}^{k}}\left[  u>v\right]
=\sum_{v=1}^{k}\left[  u>v\right] \\
&  =\sum_{v=1}^{u-1}\underbrace{\left[  u>v\right]  }%
_{\substack{=1\\\text{(since }u>v\\\text{(since }v\leq u-1<u\text{))}}%
}+\sum_{v=u}^{k}\underbrace{\left[  u>v\right]  }_{\substack{=0\\\text{(since
we don't have }u>v\\\text{(since }v\geq u\text{))}}%
}\ \ \ \ \ \ \ \ \ \ \left(  \text{since }1\leq u\leq k\right) \\
&  =\sum_{v=1}^{u-1}1+\underbrace{\sum_{v=u}^{k}0}_{=0}=\sum_{v=1}%
^{u-1}1=\left(  u-1\right)  1=u-1.
\end{align*}
This proves (\ref{pf.lem.sol.Ialbe.II.5}).}.

Now, (\ref{pf.lem.sol.Ialbe.II.4}) becomes%
\begin{align*}
\sum_{x\in I}\sum_{y\in I}\left[  x>y\right]   &  =\underbrace{\sum
_{u\in\left[  k\right]  }}_{\substack{=\sum_{u\in\left\{  1,2,\ldots
,k\right\}  }\\\text{(since }\left[  k\right]  =\left\{  1,2,\ldots,k\right\}
\text{)}}}\underbrace{\sum_{v\in\left[  k\right]  }\left[  u>v\right]
}_{\substack{=u-1\\\text{(by (\ref{pf.lem.sol.Ialbe.II.5}))}}%
}=\underbrace{\sum_{u\in\left\{  1,2,\ldots,k\right\}  }}_{=\sum_{u=1}^{k}%
}\left(  u-1\right) \\
&  =\sum_{u=1}^{k}\left(  u-1\right) \\
&  =\sum_{u=0}^{k-1}u\ \ \ \ \ \ \ \ \ \ \left(  \text{here, we have
substituted }u\text{ for }u-1\text{ in the sum}\right) \\
&  =0+1+\cdots+\left(  k-1\right)  .
\end{align*}
This proves Lemma \ref{lem.sol.Ialbe.II}.
\end{proof}
\end{verlong}

\begin{lemma}
\label{lem.sol.Ialbe.inclist2}Let $S$ be a finite set of integers. Let
$\left(  c_{1},c_{2},\ldots,c_{s}\right)  $ be a list of all elements of $S$
(with no repetitions).

Let $p_{1},p_{2},\ldots,p_{s}$ be $s$ pairwise distinct elements of $S$. Then,
there exists a $\pi\in S_{s}$ such that $\left(  p_{1},p_{2},\ldots
,p_{s}\right)  =\left(  c_{\pi\left(  1\right)  },c_{\pi\left(  2\right)
},\ldots,c_{\pi\left(  s\right)  }\right)  $.
\end{lemma}

\begin{proof}
[Proof of Lemma \ref{lem.sol.Ialbe.inclist2}.]We know that $\left(
c_{1},c_{2},\ldots,c_{s}\right)  $ is a list of all elements of $S$. Hence,
$\left\{  c_{1},c_{2},\ldots,c_{s}\right\}  =S$.

The elements $p_{1},p_{2},\ldots,p_{s}$ are pairwise distinct. In other words,
if $i$ and $j$ are two distinct elements of $\left\{  1,2,\ldots,s\right\}  $,
then%
\begin{equation}
p_{i}\neq p_{j}. \label{pf.lem.sol.Ialbe.inclist2.1}%
\end{equation}

For every $i\in\left\{  1,2,\ldots,s\right\}  $, there exists an $h\in\left\{
1,2,\ldots,s\right\}  $ such that $p_{i}=c_{h}$%
\ \ \ \ \footnote{\textit{Proof.} Let $i\in\left\{  1,2,\ldots,s\right\}  $.
Then, $p_{i}\in S$ (since $p_{1},p_{2},\ldots,p_{s}$ are $s$ elements of $S$).
Thus, $p_{i}\in S=\left\{  c_{1},c_{2},\ldots,c_{s}\right\}  $. In other
words, there exists an $h\in\left\{  1,2,\ldots,s\right\}  $ such that
$p_{i}=c_{h}$. Qed.}. Fix such an $h$, and denote it by $h_{i}$. Thus, for
every $i\in\left\{  1,2,\ldots,s\right\}  $, we have defined an $h_{i}%
\in\left\{  1,2,\ldots,s\right\}  $ such that%
\begin{equation}
p_{i}=c_{h_{i}}. \label{pf.lem.sol.Ialbe.inclist2.2}%
\end{equation}

Define a map $\varphi:\left\{  1,2,\ldots,s\right\}  \rightarrow\left\{
1,2,\ldots,s\right\}  $ by%
\[
\left(  \varphi\left(  i\right)  =h_{i}\ \ \ \ \ \ \ \ \ \ \text{for every
}i\in\left\{  1,2,\ldots,s\right\}  \right)  .
\]
Then, the map $\varphi$ is injective\footnote{\textit{Proof.} Let $i$ and $j$
be two elements of $\left\{  1,2,\ldots,s\right\}  $ such that $\varphi\left(
i\right)  =\varphi\left(  j\right)  $. We shall prove that $i=j$.
\par
The definition of $\varphi$ yields $\varphi\left(  i\right)  =h_{i}$ and
$\varphi\left(  j\right)  =h_{j}$. But (\ref{pf.lem.sol.Ialbe.inclist2.2})
yields $p_{i}=c_{h_{i}}$. Also, (\ref{pf.lem.sol.Ialbe.inclist2.2}) (applied
to $j$ instead of $i$) yields $p_{j}=c_{h_{j}}$. Now, $p_{i}=c_{h_{i}%
}=c_{h_{j}}$ (since $h_{i}=\varphi\left(  i\right)  =\varphi\left(  j\right)
=h_{j}$). Comparing this with $p_{j}=c_{h_{j}}$, we obtain $p_{i}=p_{j}$.
\par
If the elements $i$ and $j$ were distinct, then we would have $p_{i}\neq
p_{j}$ (by (\ref{pf.lem.sol.Ialbe.inclist2.1})); but this would contradict
$p_{i}=p_{j}$. Hence, the elements $i$ and $j$ cannot be distinct. In other
words, we have $i=j$.
\par
Now, forget that we fixed $i$ and $j$. We thus have proven that if $i$ and $j$
are two elements of $\left\{  1,2,\ldots,s\right\}  $ such that $\varphi
\left(  i\right)  =\varphi\left(  j\right)  $, then $i=j$. In other words, the
map $\varphi$ is injective. Qed.} and therefore
bijective\footnote{\textit{Proof.} The set $\left\{  1,2,\ldots,s\right\}  $
is finite and satisfies $\left\vert \left\{  1,2,\ldots,s\right\}  \right\vert
\geq\left\vert \left\{  1,2,\ldots,s\right\}  \right\vert $. Hence, Lemma
\ref{lem.jectivity.pigeon-inj} (applied to $U=\left\{  1,2,\ldots,s\right\}
$, $V=\left\{  1,2,\ldots,s\right\}  $ and $f=\varphi$) shows that we have the
following logical equivalence:%
\[
\left(  \varphi\text{ is injective}\right)  \ \Longleftrightarrow\ \left(
\varphi\text{ is bijective}\right)  .
\]
Hence, $\varphi$ is bijective (since $\varphi$ is injective). Qed.}. Thus,
$\varphi$ is a bijective map $\left\{  1,2,\ldots,s\right\}  \rightarrow
\left\{  1,2,\ldots,s\right\}  $. In other words, $\varphi$ is a permutation
of the set $\left\{  1,2,\ldots,s\right\}  $. In other words, $\varphi\in
S_{s}$ (since $S_{s}$ is the set of all permutations of the set $\left\{
1,2,\ldots,s\right\}  $). Furthermore, $\left(  p_{1},p_{2},\ldots
,p_{s}\right)  =\left(  c_{\varphi\left(  1\right)  },c_{\varphi\left(
2\right)  },\ldots,c_{\varphi\left(  s\right)  }\right)  $%
\ \ \ \ \footnote{\textit{Proof.} For every $i\in\left\{  1,2,\ldots
,s\right\}  $, we have $\varphi\left(  i\right)  =h_{i}$ (by the definition of
$\varphi$) and thus $c_{\varphi\left(  i\right)  }=c_{h_{i}}=p_{i}$ (by
(\ref{pf.lem.sol.Ialbe.inclist2.2})). In other words, for every $i\in\left\{
1,2,\ldots,s\right\}  $, we have $p_{i}=c_{\varphi\left(  i\right)  }$. In
other words, we have $\left(  p_{1},p_{2},\ldots,p_{s}\right)  =\left(
c_{\varphi\left(  1\right)  },c_{\varphi\left(  2\right)  },\ldots
,c_{\varphi\left(  s\right)  }\right)  $. Qed.}. Hence, $\varphi$ is an
element of $S_{s}$ and satisfies $\left(  p_{1},p_{2},\ldots,p_{s}\right)
=\left(  c_{\varphi\left(  1\right)  },c_{\varphi\left(  2\right)  }%
,\ldots,c_{\varphi\left(  s\right)  }\right)  $. Thus, there exists a $\pi\in
S_{s}$ such that $\left(  p_{1},p_{2},\ldots,p_{s}\right)  =\left(
c_{\pi\left(  1\right)  },c_{\pi\left(  2\right)  },\ldots,c_{\pi\left(
s\right)  }\right)  $ (namely, $\pi=\varphi$). This proves Lemma
\ref{lem.sol.Ialbe.inclist2}.
\end{proof}

\begin{lemma}
\label{lem.sol.Ialbe.InotI}Let $n\in\mathbb{N}$. Let $I$ be a subset of
$\left[  n\right]  $. Let $k=\left\vert I\right\vert $.

Let $\left(  a_{1},a_{2},\ldots,a_{k}\right)  $ be a list of all elements of
$I$ (with no repetitions). Let $\left(  b_{1},b_{2},\ldots,b_{n-k}\right)  $
be a list of all elements of $\left[  n\right]  \setminus I$ (with no repetitions).

\textbf{(a)} There exists a unique $\sigma\in S_{n}$ satisfying%
\[
\left(  \sigma\left(  1\right)  ,\sigma\left(  2\right)  ,\ldots,\sigma\left(
n\right)  \right)  =\left(  a_{1},a_{2},\ldots,a_{k},b_{1},b_{2}%
,\ldots,b_{n-k}\right)  .
\]

\textbf{(b)} Let $\sum I$ denote the sum of all elements of $I$. (Thus, $\sum
I=\sum_{i\in I}i$.) Then,%
\[
\sum_{\left(  i,j\right)  \in\left[  k\right]  \times\left[  n-k\right]
}\left[  a_{i}>b_{j}\right]  =\sum I-\left(  1+2+\cdots+k\right)  .
\]

\end{lemma}

\begin{vershort}
\begin{proof}
[Proof of Lemma \ref{lem.sol.Ialbe.InotI}.]\textbf{(a)} The $k$ elements
$a_{1},a_{2},\ldots,a_{k}$ belong to $I$ (since $\left(  a_{1},a_{2}%
,\ldots,a_{k}\right)  $ is a list of all elements of $I$), and thus belong to
$\left[  n\right]  $ as well (since $I\subseteq\left[  n\right]  $). Also, the
$k$ elements $a_{1},a_{2},\ldots,a_{k}$ are pairwise distinct (since $\left(
a_{1},a_{2},\ldots,a_{k}\right)  $ is a list with no repetitions).

The $n-k$ elements $b_{1},b_{2},\ldots,b_{n-k}$ belong to $\left[  n\right]
\setminus I$ (since $\left(  b_{1},b_{2},\ldots,b_{n-k}\right)  $ is a list of
all elements of $\left[  n\right]  \setminus I$), and thus belong to $\left[
n\right]  $ as well (since $\left[  n\right]  \setminus I\subseteq\left[
n\right]  $). Also, the $n-k$ elements $b_{1},b_{2},\ldots,b_{n-k}$ are
pairwise distinct (since $\left(  b_{1},b_{2},\ldots,b_{n-k}\right)  $ is a
list with no repetitions).

The $n$ elements $a_{1},a_{2},\ldots,a_{k},b_{1},b_{2},\ldots,b_{n-k}$ are
pairwise distinct\footnote{\textit{Proof.} Assume the contrary. Thus, two of
the $n$ elements $a_{1},a_{2},\ldots,a_{k},b_{1},b_{2},\ldots,b_{n-k}$ are
equal. These two elements are either two of the elements $a_{1},a_{2}%
,\ldots,a_{k}$, or two of the elements $b_{1},b_{2},\ldots,b_{n-k}$, or one of
the former and one of the latter (in some order). Hence, we are in one of the
following three cases:
\par
\textit{Case 1:} Two of the $k$ elements $a_{1},a_{2},\ldots,a_{k}$ are equal.
\par
\textit{Case 2:} Two of the $n-k$ elements $b_{1},b_{2},\ldots,b_{n-k}$ are
equal.
\par
\textit{Case 3:} One of the $k$ elements $a_{1},a_{2},\ldots,a_{k}$ is equal
to one of the elements $b_{1},b_{2},\ldots,b_{n-k}$.
\par
But Case 1 is impossible (since the $k$ elements $a_{1},a_{2},\ldots,a_{k}$
are pairwise distinct), and Case 2 is also impossible (since the $n-k$
elements $b_{1},b_{2},\ldots,b_{n-k}$ are pairwise distinct). Hence, we must
be in Case 3. In other words, one of the $k$ elements $a_{1},a_{2}%
,\ldots,a_{k}$ is equal to one of the elements $b_{1},b_{2},\ldots,b_{n-k}$.
In other words, we must have $a_{i}=b_{j}$ for some $i\in\left\{
1,2,\ldots,k\right\}  $ and some $j\in\left\{  1,2,\ldots,n-k\right\}  $.
Consider these $i$ and $j$. We have $a_{i}\in I$ (since the $k$ elements
$a_{1},a_{2},\ldots,a_{k}$ belong to $I$). But $b_{j}\in\left[  n\right]
\setminus I$ (since the $n-k$ elements $b_{1},b_{2},\ldots,b_{n-k}$ belong to
$\left[  n\right]  \setminus I$) and thus $b_{j}\notin I$. Hence, $a_{i}%
=b_{j}\notin I$; but this contradicts $a_{i}\in I$. This contradiction proves
that our assumption was wrong; qed.} and belong to $\left[  n\right]
$\ \ \ \ \footnote{This is because the $k$ elements $a_{1},a_{2},\ldots,a_{k}$
belong to $\left[  n\right]  $, and because the $n-k$ elements $b_{1}%
,b_{2},\ldots,b_{n-k}$ belong to $\left[  n\right]  $.}. Furthermore, $\left(
1,2,\ldots,n\right)  $ is a list of all elements of $\left[  n\right]  $ (with
no repetition). Hence, Lemma \ref{lem.sol.Ialbe.inclist2} (applied to
$S=\left[  n\right]  $, $s=n$, $\left(  c_{1},c_{2},\ldots,c_{s}\right)
=\left(  1,2,\ldots,n\right)  $ and $\left(  p_{1},p_{2},\ldots,p_{s}\right)
=\left(  a_{1},a_{2},\ldots,a_{k},b_{1},b_{2},\ldots,b_{n-k}\right)  $) yields
that there exists a $\pi\in S_{n}$ such that
\[
\left(  a_{1},a_{2},\ldots,a_{k},b_{1},b_{2},\ldots,b_{n-k}\right)  =\left(
\pi\left(  1\right)  ,\pi\left(  2\right)  ,\ldots,\pi\left(  n\right)
\right)  .
\]
Consider this $\pi$.

We have $\pi\in S_{n}$ and $\left(  \pi\left(  1\right)  ,\pi\left(  2\right)
,\ldots,\pi\left(  n\right)  \right)  =\left(  a_{1},a_{2},\ldots,a_{k}%
,b_{1},b_{2},\ldots,b_{n-k}\right)  $.

But our goal is to show that there exists a unique $\sigma\in S_{n}$
satisfying%
\begin{equation}
\left(  \sigma\left(  1\right)  ,\sigma\left(  2\right)  ,\ldots,\sigma\left(
n\right)  \right)  =\left(  a_{1},a_{2},\ldots,a_{k},b_{1},b_{2}%
,\ldots,b_{n-k}\right)  . \label{pf.lem.sol.Ialbe.InotI.short.a.1}%
\end{equation}
We already know that there exists \textbf{at least one} such $\sigma$ (namely,
$\sigma=\pi$). Hence, it remains to show that there exists \textbf{at most
one} such $\sigma$. In other words, it remains to show that if $\sigma_{1}$
and $\sigma_{2}$ are two elements $\sigma\in S_{n}$ satisfying
(\ref{pf.lem.sol.Ialbe.InotI.short.a.1}), then $\sigma_{1}=\sigma_{2}$.

So let $\sigma_{1}$ and $\sigma_{2}$ be two elements $\sigma\in S_{n}$
satisfying (\ref{pf.lem.sol.Ialbe.InotI.short.a.1}). We must prove that
$\sigma_{1}=\sigma_{2}$.

We know that $\sigma_{1}$ is an element $\sigma\in S_{n}$ satisfying
(\ref{pf.lem.sol.Ialbe.InotI.short.a.1}). In other words, $\sigma_{1}$ is an
element of $S_{n}$ and satisfies%
\begin{equation}
\left(  \sigma_{1}\left(  1\right)  ,\sigma_{1}\left(  2\right)
,\ldots,\sigma_{1}\left(  n\right)  \right)  =\left(  a_{1},a_{2},\ldots
,a_{k},b_{1},b_{2},\ldots,b_{n-k}\right)  .
\label{pf.lem.sol.Ialbe.InotI.short.a.3a}%
\end{equation}
Similarly, $\sigma_{2}$ is an element of $S_{n}$ and satisfies%
\begin{equation}
\left(  \sigma_{2}\left(  1\right)  ,\sigma_{2}\left(  2\right)
,\ldots,\sigma_{2}\left(  n\right)  \right)  =\left(  a_{1},a_{2},\ldots
,a_{k},b_{1},b_{2},\ldots,b_{n-k}\right)  .
\label{pf.lem.sol.Ialbe.InotI.short.a.3b}%
\end{equation}
Comparing (\ref{pf.lem.sol.Ialbe.InotI.short.a.3a}) with
(\ref{pf.lem.sol.Ialbe.InotI.short.a.3b}), we obtain%
\[
\left(  \sigma_{1}\left(  1\right)  ,\sigma_{1}\left(  2\right)
,\ldots,\sigma_{1}\left(  n\right)  \right)  =\left(  \sigma_{2}\left(
1\right)  ,\sigma_{2}\left(  2\right)  ,\ldots,\sigma_{2}\left(  n\right)
\right)  .
\]
In other words, $\sigma_{1}\left(  i\right)  =\sigma_{2}\left(  i\right)  $
for every $i\in\left[  n\right]  $. Since both $\sigma_{1}$ and $\sigma_{2}$
are maps $\left[  n\right]  \rightarrow\left[  n\right]  $, this shows that
$\sigma_{1}=\sigma_{2}$. This is precisely what we wanted to show.

Thus, we have proven that there exists \textbf{at most one} $\sigma\in S_{n}$
satisfying (\ref{pf.lem.sol.Ialbe.InotI.short.a.1}). As explained above, this
completes the proof of Lemma \ref{lem.sol.Ialbe.InotI} \textbf{(a)}.

\textbf{(b)} Clearly,
\begin{equation}
\sum I=\sum_{i\in I}i=\sum_{x\in I}x
\label{pf.lem.sol.Ialbe.InotI.short.b.triv}%
\end{equation}
(here, we have renamed the summation index $i$ as $x$).

We know that $\left(  a_{1},a_{2},\ldots,a_{k}\right)  $ is a list of all
elements of $I$ (with no repetitions). Therefore, Lemma
\ref{lem.sol.Ialbe.inclist} \textbf{(a)} (applied to $I$, $k$ and $\left(
a_{1},a_{2},\ldots,a_{k}\right)  $ instead of $S$, $s$ and $\left(
c_{1},c_{2},\ldots,c_{s}\right)  $) shows that the map $\left[  k\right]
\rightarrow I,\ h\mapsto a_{h}$ is well-defined and a bijection.

We know that $\left(  b_{1},b_{2},\ldots,b_{n-k}\right)  $ is a list of all
elements of $\left[  n\right]  \setminus I$ (with no repetitions). Therefore,
Lemma \ref{lem.sol.Ialbe.inclist} \textbf{(a)} (applied to $\left[  n\right]
\setminus I$, $n-k$ and $\left(  b_{1},b_{2},\ldots,b_{n-k}\right)  $ instead
of $S$, $s$ and $\left(  c_{1},c_{2},\ldots,c_{s}\right)  $) shows that the
map $\left[  n-k\right]  \rightarrow\left[  n\right]  \setminus I,\ h\mapsto
b_{h}$ is well-defined and a bijection.

Now,%
\begin{align}
\underbrace{\sum_{\left(  i,j\right)  \in\left[  k\right]  \times\left[
n-k\right]  }}_{=\sum_{i\in\left[  k\right]  }\sum_{j\in\left[  n-k\right]  }%
}\left[  a_{i}>b_{j}\right]   &  =\sum_{i\in\left[  k\right]  }%
\underbrace{\sum_{j\in\left[  n-k\right]  }\left[  a_{i}>b_{j}\right]
}_{\substack{=\sum_{y\in\left[  n\right]  \setminus I}\left[  a_{i}>y\right]
\\\text{(here, we have substituted }y\text{ for }b_{j}\\\text{in the sum,
since the}\\\text{map }\left[  n-k\right]  \rightarrow\left[  n\right]
\setminus I,\ h\mapsto b_{h}\text{ is a bijection)}}}=\sum_{i\in\left[
k\right]  }\sum_{y\in\left[  n\right]  \setminus I}\left[  a_{i}>y\right]
\nonumber\\
&  =\sum_{x\in I}\sum_{y\in\left[  n\right]  \setminus I}\left[  x>y\right]
\label{pf.lem.sol.Ialbe.InotI.short.b.1}%
\end{align}
(here, we have substituted $x$ for $a_{i}$ in the outer sum, since the map
$\left[  k\right]  \rightarrow I,\ h\mapsto a_{h}$ is a bijection). But every
$x\in I$ satisfies%
\begin{equation}
\sum_{y\in\left[  n\right]  \setminus I}\left[  x>y\right]  =x-1-\sum_{y\in
I}\left[  x>y\right]  \label{pf.lem.sol.Ialbe.InotI.short.b.sum1}%
\end{equation}
\footnote{\textit{Proof of (\ref{pf.lem.sol.Ialbe.InotI.short.b.sum1}):} Let
$x\in I$. Then, $x\in I\subseteq\left[  n\right]  $, so that $1\leq x\leq n$.
\par
Every $y\in\left[  n\right]  $ satisfies either $y\in I$ or $y\notin I$ (but
not both). Hence,%
\[
\sum_{y\in\left[  n\right]  }\left[  x>y\right]  =\underbrace{\sum
_{\substack{y\in\left[  n\right]  ;\\y\in I}}}_{\substack{=\sum_{y\in
I}\\\text{(since }I\subseteq\left[  n\right]  \text{)}}}\left[  x>y\right]
+\underbrace{\sum_{\substack{y\in\left[  n\right]  ;\\y\notin I}}}%
_{=\sum_{y\in\left[  n\right]  \setminus I}}\left[  x>y\right]  =\sum_{y\in
I}\left[  x>y\right]  +\sum_{y\in\left[  n\right]  \setminus I}\left[
x>y\right]  .
\]
Comparing this with%
\begin{align*}
\underbrace{\sum_{y\in\left[  n\right]  }}_{=\sum_{y=1}^{n}}\left[
x>y\right]   &  =\sum_{y=1}^{n}\left[  x>y\right] \\
&  =\sum_{y=1}^{x-1}\underbrace{\left[  x>y\right]  }%
_{\substack{=1\\\text{(since }x>y\\\text{(since }y\leq x-1<x\text{))}}%
}+\sum_{y=x}^{n}\underbrace{\left[  x>y\right]  }_{\substack{=0\\\text{(since
we don't have }x>y\\\text{(since }y\geq x\text{))}}%
}\ \ \ \ \ \ \ \ \ \ \left(  \text{since }1\leq x\leq n\right) \\
&  =\sum_{y=1}^{x-1}1+\underbrace{\sum_{y=x}^{n}0}_{=0}=\sum_{y=1}%
^{x-1}1=\left(  x-1\right)  1=x-1,
\end{align*}
we obtain%
\[
\sum_{y\in I}\left[  x>y\right]  +\sum_{y\in\left[  n\right]  \setminus
I}\left[  x>y\right]  =x-1.
\]
Subtracting $\sum_{y\in I}\left[  x>y\right]  $ from both sides of this
equality, we obtain $\sum_{y\in\left[  n\right]  \setminus I}\left[
x>y\right]  =x-1-\sum_{y\in I}\left[  x>y\right]  $. This proves
(\ref{pf.lem.sol.Ialbe.InotI.short.b.sum1}).}.

Now, (\ref{pf.lem.sol.Ialbe.InotI.short.b.1}) becomes%
\begin{align*}
&  \sum_{\left(  i,j\right)  \in\left[  k\right]  \times\left[  n-k\right]
}\left[  a_{i}>b_{j}\right] \\
&  =\sum_{x\in I}\underbrace{\sum_{y\in\left[  n\right]  \setminus I}\left[
x>y\right]  }_{\substack{=x-1-\sum_{y\in I}\left[  x>y\right]  \\\text{(by
(\ref{pf.lem.sol.Ialbe.InotI.short.b.sum1}))}}}=\sum_{x\in I}\left(
x-1-\sum_{y\in I}\left[  x>y\right]  \right) \\
&  =\underbrace{\sum_{x\in I}x}_{\substack{=\sum I\\\text{(by
(\ref{pf.lem.sol.Ialbe.InotI.short.b.triv}))}}}-\underbrace{\sum_{x\in I}%
1}_{=\left\vert I\right\vert \cdot1=\left\vert I\right\vert =k}%
-\underbrace{\sum_{x\in I}\sum_{y\in I}\left[  x>y\right]  }%
_{\substack{=0+1+\cdots+\left(  k-1\right)  \\\text{(by Lemma
\ref{lem.sol.Ialbe.II})}}}=\sum I-k-\left(  0+1+\cdots+\left(  k-1\right)
\right) \\
&  =\sum I-\underbrace{\left(  \left(  0+1+\cdots+\left(  k-1\right)  \right)
+k\right)  }_{=0+1+\cdots+k=1+2+\cdots+k}=\sum I-\left(  1+2+\cdots+k\right)
.
\end{align*}
This proves Lemma \ref{lem.sol.Ialbe.InotI} \textbf{(b)}.
\end{proof}
\end{vershort}

\begin{verlong}
\begin{proof}
[Proof of Lemma \ref{lem.sol.Ialbe.InotI}.]\textbf{(a)} The list $\left(
a_{1},a_{2},\ldots,a_{k}\right)  $ is a list of all elements of $I$, and thus
belongs to $I^{k}$. Hence, $\left(  a_{1},a_{2},\ldots,a_{k}\right)  \in
I^{k}\subseteq\left[  n\right]  ^{k}$ (since $I\subseteq\left[  n\right]  $).

The list $\left(  b_{1},b_{2},\ldots,b_{n-k}\right)  $ is a list of all
elements of $\left[  n\right]  \setminus I$, and thus belongs to $\left(
\left[  n\right]  \setminus I\right)  ^{n-k}$. Hence, $\left(  b_{1}%
,b_{2},\ldots,b_{n-k}\right)  \in\left(  \left[  n\right]  \setminus I\right)
^{n-k}\subseteq\left[  n\right]  ^{n-k}$ (since $\left[  n\right]  \setminus
I\subseteq\left[  n\right]  $).

From $\left(  a_{1},a_{2},\ldots,a_{k}\right)  \in\left[  n\right]  ^{k}$ and
$\left(  b_{1},b_{2},\ldots,b_{n-k}\right)  \in\left[  n\right]  ^{n-k}$, we
obtain
\[
\left(  a_{1},a_{2},\ldots,a_{k},b_{1},b_{2},\ldots,b_{n-k}\right)  \in\left[
n\right]  ^{k+\left(  n-k\right)  }=\left[  n\right]  ^{n}%
\ \ \ \ \ \ \ \ \ \ \left(  \text{since }k+\left(  n-k\right)  =n\right)  .
\]
In other words, $\left(  a_{1},a_{2},\ldots,a_{k},b_{1},b_{2},\ldots
,b_{n-k}\right)  $ is an $n$-tuple of elements of $\left[  n\right]  $. Denote
this $n$-tuple by $\left(  c_{1},c_{2},\ldots,c_{n}\right)  $. Thus,%
\[
\left(  c_{1},c_{2},\ldots,c_{n}\right)  =\left(  a_{1},a_{2},\ldots
,a_{k},b_{1},b_{2},\ldots,b_{n-k}\right)  .
\]
In other words,%
\[
c_{h}=%
\begin{cases}
a_{h}, & \text{if }h\leq k;\\
b_{h-k}, & \text{if }h>k
\end{cases}
\ \ \ \ \ \ \ \ \ \ \text{for every }h\in\left\{  1,2,\ldots,n\right\}  .
\]
In other words,%
\begin{equation}
c_{h}=%
\begin{cases}
a_{h}, & \text{if }h\leq k;\\
b_{h-k}, & \text{if }h>k
\end{cases}
\ \ \ \ \ \ \ \ \ \ \text{for every }h\in\left[  n\right]
\label{pf.lem.sol.Ialbe.InotI.a.ch=}%
\end{equation}
(since $\left\{  1,2,\ldots,n\right\}  =\left[  n\right]  $). Notice that%
\begin{equation}
c_{h}\notin I\ \ \ \ \ \ \ \ \ \ \text{for every }h\in\left[  n\right]  \text{
satisfying }h>k \label{pf.lem.sol.Ialbe.InotI.a.notinI}%
\end{equation}
\footnote{\textit{Proof of (\ref{pf.lem.sol.Ialbe.InotI.a.notinI}):} Let
$h\in\left[  n\right]  $ be such that $h>k$. Then,%
\begin{align*}
c_{h}  &  =%
\begin{cases}
a_{h}, & \text{if }h\leq k;\\
b_{h-k}, & \text{if }h>k
\end{cases}
\ \ \ \ \ \ \ \ \ \ \left(  \text{by (\ref{pf.lem.sol.Ialbe.InotI.a.ch=}%
)}\right) \\
&  =b_{h-k}\ \ \ \ \ \ \ \ \ \ \left(  \text{since }h>k\right) \\
&  \in\left[  n\right]  \setminus I\ \ \ \ \ \ \ \ \ \ \left(  \text{since
}\left(  b_{1},b_{2},\ldots,b_{n-k}\right)  \in\left(  \left[  n\right]
\setminus I\right)  ^{n-k}\right)  .
\end{align*}
In other words, $c_{h}\in\left[  n\right]  $ and $c_{h}\notin I$. This proves
(\ref{pf.lem.sol.Ialbe.InotI.a.notinI}).}.

It is well-known that if $U$ is a set, and if $V$ is a subset of $U$, then
$V\cup\left(  U\setminus V\right)  =U$. Applying this to $U=\left[  n\right]
$ and $V=I$, we obtain $I\cup\left(  \left[  n\right]  \setminus I\right)
=\left[  n\right]  $ (since $I$ is a subset of $\left[  n\right]  $).

Now, define a map $f:\left[  n\right]  \rightarrow\left[  n\right]  $ by%
\[
\left(  f\left(  h\right)  =c_{h}\ \ \ \ \ \ \ \ \ \ \text{for every }%
h\in\left[  n\right]  \right)  .
\]
(This is well-defined, because $c_{h}$ is a well-defined element of $\left[
n\right]  $ for every $h\in\left[  n\right]  $\ \ \ \ \footnote{Indeed, this
follows from (\ref{pf.lem.sol.Ialbe.InotI.a.ch=}).}.) We have%
\begin{align*}
\left(  f\left(  1\right)  ,f\left(  2\right)  ,\ldots,f\left(  n\right)
\right)   &  =\left(  c_{1},c_{2},\ldots,c_{n}\right)
\ \ \ \ \ \ \ \ \ \ \left(  \text{since }f\left(  h\right)  =c_{h}\text{ for
every }h\in\left[  n\right]  \right) \\
&  =\left(  a_{1},a_{2},\ldots,a_{k},b_{1},b_{2},\ldots,b_{n-k}\right)
\end{align*}
and thus%
\begin{align*}
\left\{  f\left(  1\right)  ,f\left(  2\right)  ,\ldots,f\left(  n\right)
\right\}   &  =\left\{  a_{1},a_{2},\ldots,a_{k},b_{1},b_{2},\ldots
,b_{n-k}\right\} \\
&  =\underbrace{\left\{  a_{1},a_{2},\ldots,a_{k}\right\}  }%
_{\substack{=I\\\text{(since }\left(  a_{1},a_{2},\ldots,a_{k}\right)  \text{
is a list}\\\text{of all elements of }I\text{)}}}\cup\underbrace{\left\{
b_{1},b_{2},\ldots,b_{n-k}\right\}  }_{\substack{=\left[  n\right]  \setminus
I\\\text{(since }\left(  b_{1},b_{2},\ldots,b_{n-k}\right)  \text{ is a
list}\\\text{of all elements of }\left[  n\right]  \setminus I\text{)}}}\\
&  =I\cup\left(  \left[  n\right]  \setminus I\right)  =\left[  n\right]  .
\end{align*}
Now,%
\[
f\left(  \underbrace{\left[  n\right]  }_{=\left\{  1,2,\ldots,n\right\}
}\right)  =f\left(  \left\{  1,2,\ldots,n\right\}  \right)  =\left\{  f\left(
1\right)  ,f\left(  2\right)  ,\ldots,f\left(  n\right)  \right\}  =\left[
n\right]  ,
\]
so that $\left[  n\right]  \subseteq\left[  n\right]  =f\left(  \left[
n\right]  \right)  $. In other words, the map $f:\left[  n\right]
\rightarrow\left[  n\right]  $ is surjective. From this, it is easy to
conclude that this map $f$ is injective\footnote{\textit{First proof.} The set
$\left[  n\right]  $ is finite and satisfies $\left\vert \left[  n\right]
\right\vert \leq\left\vert \left[  n\right]  \right\vert $. Hence, Lemma
\ref{lem.jectivity.pigeon-surj} (applied to $\left[  n\right]  $ and $\left[
n\right]  $ instead of $U$ and $V$) shows that we have the following logical
equivalence:%
\[
\left(  f\text{ is surjective}\right)  \ \Longleftrightarrow\ \left(  f\text{
is bijective}\right)  .
\]
Thus, $f$ is bijective (since $f$ is surjective). Hence, $f$ is injective.
Qed.
\par
\textit{Second proof.} Let us give a different proof of the injectivity of
$f$, which does not use the finiteness of $\left[  n\right]  $.
\par
Let $g$ and $h$ be two elements of $\left[  n\right]  $ such that $f\left(
g\right)  =f\left(  h\right)  $. We shall prove that $g=h$.
\par
We can WLOG assume that $g\leq h$ (since otherwise, we can simply swap $g$
with $h$). Assume this.
\par
The definition of $f$ yields $f\left(  g\right)  =c_{g}$ and $f\left(
h\right)  =c_{h}$. Thus, $c_{g}=f\left(  g\right)  =f\left(  h\right)  =c_{h}%
$.
\par
For the sake of contradiction, let us assume that $g\neq h$. Combining this
with $g\leq h$, we obtain $g<h$. We are in one of the following two cases:
\par
\textit{Case 1:} We have $g\leq k$.
\par
\textit{Case 2:} We have $g>k$.
\par
Let us first consider Case 1. In this case, we have $g\leq k$. Thus,
(\ref{pf.lem.sol.Ialbe.InotI.a.ch=}) (applied to $g$ instead of $h$) yields
$c_{g}=%
\begin{cases}
a_{g}, & \text{if }g\leq k;\\
b_{g-k}, & \text{if }g>k
\end{cases}
=a_{g}$ (since $g\leq k$). Thus, $c_{g}=a_{g}\in I$ (since $\left(
a_{1},a_{2},\ldots,a_{k}\right)  \in I^{k}$). If we had $h>k$, then we would
have%
\[
c_{g}=c_{h}\notin I\ \ \ \ \ \ \ \ \ \ \left(  \text{by
(\ref{pf.lem.sol.Ialbe.InotI.a.notinI})}\right)  ,
\]
which would contradict $c_{g}\in I$. Thus, we cannot have $h>k$. Hence, we
must have $h\leq k$. Thus, (\ref{pf.lem.sol.Ialbe.InotI.a.ch=}) yields $c_{h}=%
\begin{cases}
a_{h}, & \text{if }h\leq k;\\
b_{h-k}, & \text{if }h>k
\end{cases}
=a_{h}$ (since $h\leq k$). Now, $c_{g}=a_{g}$, so that $a_{g}=c_{g}%
=c_{h}=a_{h}$.
\par
From $g\in\left[  n\right]  =\left\{  1,2,\ldots,n\right\}  $, we obtain
$g\geq1$. Combining this with $g\leq k$, we obtain $g\in\left\{
1,2,\ldots,k\right\}  $.
\par
From $h\in\left[  n\right]  =\left\{  1,2,\ldots,n\right\}  $, we obtain
$h\geq1$. Combining this with $h\leq k$, we obtain $h\in\left\{
1,2,\ldots,k\right\}  $.
\par
Now, $g$ and $h$ are two elements of $\left\{  1,2,\ldots,k\right\}  $ (since
$g\in\left\{  1,2,\ldots,k\right\}  $ and $h\in\left\{  1,2,\ldots,k\right\}
$). These two elements are distinct (since $g<h$).
\par
But $\left(  a_{1},a_{2},\ldots,a_{k}\right)  $ is a list of all elements of
$I$ (with no repetitions). Thus, the list $\left(  a_{1},a_{2},\ldots
,a_{k}\right)  $ contains no repetitions. In other words, the elements
$a_{1},a_{2},\ldots,a_{k}$ are pairwise distinct. In other words, if $u$ and
$v$ are two distinct elements of $\left\{  1,2,\ldots,k\right\}  $, then
$a_{u}\neq a_{v}$. Applying this to $u=g$ and $v=h$, we obtain $a_{g}\neq
a_{h}$. This contradicts $a_{g}=a_{h}$. Hence, we have obtained a
contradiction in Case 1.
\par
Let us now consider Case 2. In this case, we have $g>k$. Hence, $g\geq k+1$
(since $g$ and $k$ are integers). But (\ref{pf.lem.sol.Ialbe.InotI.a.ch=})
(applied to $g$ instead of $h$) yields $c_{g}=%
\begin{cases}
a_{g}, & \text{if }g\leq k;\\
b_{g-k}, & \text{if }g>k
\end{cases}
=b_{g-k}$ (since $g>k$). On the other hand, $g\leq h$, so that $h\geq g\geq
k+1>k$. Thus, (\ref{pf.lem.sol.Ialbe.InotI.a.ch=}) yields $c_{h}=%
\begin{cases}
a_{h}, & \text{if }h\leq k;\\
b_{h-k}, & \text{if }h>k
\end{cases}
=b_{h-k}$ (since $h>k$). Now, $c_{g}=b_{g-k}$, so that $b_{g-k}=c_{g}%
=c_{h}=b_{h-k}$.
\par
From $g\in\left[  n\right]  =\left\{  1,2,\ldots,n\right\}  $, we obtain
$g\leq n$. Combining this with $g\geq k+1$, we obtain $g\in\left\{
k+1,k+2,\ldots,n\right\}  $. Hence, $g-k\in\left\{  1,2,\ldots,n-k\right\}  $.
\par
From $h\in\left[  n\right]  =\left\{  1,2,\ldots,n\right\}  $, we obtain
$h\leq n$. Combining this with $h\geq k+1$, we obtain $h\in\left\{
k+1,k+2,\ldots,n\right\}  $. Hence, $h-k\in\left\{  1,2,\ldots,n-k\right\}  $.
\par
Now, $g-k$ and $h-k$ are two elements of $\left\{  1,2,\ldots,n-k\right\}  $
(since $g-k\in\left\{  1,2,\ldots,n-k\right\}  $ and $h-k\in\left\{
1,2,\ldots,n-k\right\}  $). These two elements are distinct (since
$\underbrace{g}_{<h}-k<h-k$).
\par
But $\left(  b_{1},b_{2},\ldots,b_{n-k}\right)  $ is a list of all elements of
$\left[  n\right]  \setminus I$ (with no repetitions). Thus, the list $\left(
b_{1},b_{2},\ldots,b_{n-k}\right)  $ contains no repetitions. In other words,
the elements $b_{1},b_{2},\ldots,b_{n-k}$ are pairwise distinct. In other
words, if $u$ and $v$ are two distinct elements of $\left\{  1,2,\ldots
,n-k\right\}  $, then $b_{u}\neq b_{v}$. Applying this to $u=g-k$ and $v=h-k$,
we obtain $b_{g-k}\neq b_{h-k}$. This contradicts $b_{g-k}=b_{h-k}$. Hence, we
have obtained a contradiction in Case 2.
\par
Now, we have obtained a contradiction in each of the two Cases 1 and 2. Thus,
we always have a contradiction (since the two Cases 1 and 2 cover all
possibilities). This contradiction proves that our assumption (that $g\neq h$)
was wrong. Hence, we cannot have $g\neq h$. We thus have $g=h$.
\par
Now, forget that we fixed $g$ and $h$. We thus have proven that if $g$ and $h$
are two elements of $\left[  n\right]  $ such that $f\left(  g\right)
=f\left(  h\right)  $, then $g=h$. In other words, the map $f$ is injective.
Qed.}. Hence, the map $f$ is bijective (since $f$ is injective and
surjective). Hence, $f$ is a bijective map $\left[  n\right]  \rightarrow
\left[  n\right]  $. In other words, $f$ is a permutation of $\left[
n\right]  $.

But $S_{n}$ is the set of all permutations of $\left\{  1,2,\ldots,n\right\}
$ (by the definition of $S_{n}$). In other words, $S_{n}$ is the set of all
permutations of $\left[  n\right]  $ (since $\left[  n\right]  =\left\{
1,2,\ldots,n\right\}  $). Now, $f$ is a permutation of $\left[  n\right]  $.
In other words, $f$ is an element of $S_{n}$ (since $S_{n}$ is the set of all
permutations of $\left[  n\right]  $). In other words, $f\in S_{n}$. Hence,
$f$ is a $\sigma\in S_{n}$ satisfying%
\begin{equation}
\left(  \sigma\left(  1\right)  ,\sigma\left(  2\right)  ,\ldots,\sigma\left(
n\right)  \right)  =\left(  a_{1},a_{2},\ldots,a_{k},b_{1},b_{2}%
,\ldots,b_{n-k}\right)  \label{pf.lem.sol.Ialbe.InotI.a.silis}%
\end{equation}
(since $f$ is an element of $S_{n}$ and satisfies \newline$\left(  f\left(
1\right)  ,f\left(  2\right)  ,\ldots,f\left(  n\right)  \right)  =\left(
a_{1},a_{2},\ldots,a_{k},b_{1},b_{2},\ldots,b_{n-k}\right)  $). Thus, there
exists \textbf{at least one} $\sigma\in S_{n}$ satisfying
(\ref{pf.lem.sol.Ialbe.InotI.a.silis}).

On the other hand, if $\sigma_{1}$ and $\sigma_{2}$ are two elements
$\sigma\in S_{n}$ satisfying (\ref{pf.lem.sol.Ialbe.InotI.a.silis}), then
$\sigma_{1}=\sigma_{2}$\ \ \ \ \footnote{\textit{Proof.} Let $\sigma_{1}$ and
$\sigma_{2}$ be two elements $\sigma\in S_{n}$ satisfying
(\ref{pf.lem.sol.Ialbe.InotI.a.silis}). We shall show that $\sigma_{1}%
=\sigma_{2}$.
\par
We know that $\sigma_{1}$ is an element $\sigma\in S_{n}$ satisfying
(\ref{pf.lem.sol.Ialbe.InotI.a.silis}). In other words, $\sigma_{1}$ is an
element of $S_{n}$ and satisfies%
\[
\left(  \sigma_{1}\left(  1\right)  ,\sigma_{1}\left(  2\right)
,\ldots,\sigma_{1}\left(  n\right)  \right)  =\left(  a_{1},a_{2},\ldots
,a_{k},b_{1},b_{2},\ldots,b_{n-k}\right)  .
\]
\par
We know that $\sigma_{2}$ is an element $\sigma\in S_{n}$ satisfying
(\ref{pf.lem.sol.Ialbe.InotI.a.silis}). In other words, $\sigma_{2}$ is an
element of $S_{n}$ and satisfies%
\[
\left(  \sigma_{2}\left(  1\right)  ,\sigma_{2}\left(  2\right)
,\ldots,\sigma_{2}\left(  n\right)  \right)  =\left(  a_{1},a_{2},\ldots
,a_{k},b_{1},b_{2},\ldots,b_{n-k}\right)  .
\]
\par
Now,
\[
\left(  \sigma_{1}\left(  1\right)  ,\sigma_{1}\left(  2\right)
,\ldots,\sigma_{1}\left(  n\right)  \right)  =\left(  a_{1},a_{2},\ldots
,a_{k},b_{1},b_{2},\ldots,b_{n-k}\right)  =\left(  \sigma_{2}\left(  1\right)
,\sigma_{2}\left(  2\right)  ,\ldots,\sigma_{2}\left(  n\right)  \right)  .
\]
In other words, $\sigma_{1}\left(  i\right)  =\sigma_{2}\left(  i\right)  $
for every $i\in\left\{  1,2,\ldots,n\right\}  $. In other words, $\sigma
_{1}\left(  i\right)  =\sigma_{2}\left(  i\right)  $ for every $i\in\left[
n\right]  $ (since $\left\{  1,2,\ldots,n\right\}  =\left[  n\right]  $).
\par
Now, $\sigma_{1}$ is an element of $S_{n}$. In other words, $\sigma_{1}$ is a
permutation of $\left[  n\right]  $ (since $S_{n}$ is the set of all
permutations of $\left[  n\right]  $). Hence, $\sigma_{1}$ is a map $\left[
n\right]  \rightarrow\left[  n\right]  $. The same argument (applied to
$\sigma_{2}$ instead of $\sigma_{1}$) shows that $\sigma_{2}$ is a map
$\left[  n\right]  \rightarrow\left[  n\right]  $. Hence, $\sigma_{1}$ and
$\sigma_{2}$ are maps $\left[  n\right]  \rightarrow\left[  n\right]  $. Thus,
$\sigma_{1}=\sigma_{2}$ (because $\sigma_{1}\left(  i\right)  =\sigma
_{2}\left(  i\right)  $ for every $i\in\left[  n\right]  $). Qed.}. In other
words, there exists \textbf{at most one} $\sigma\in S_{n}$ satisfying
(\ref{pf.lem.sol.Ialbe.InotI.a.silis}).

We have now proven that there exists \textbf{at least one} $\sigma\in S_{n}$
satisfying (\ref{pf.lem.sol.Ialbe.InotI.a.silis}), and that there exists
\textbf{at most one} such $\sigma$. Hence, there exists a \textbf{unique}
$\sigma\in S_{n}$ satisfying (\ref{pf.lem.sol.Ialbe.InotI.a.silis}). This
proves Lemma \ref{lem.sol.Ialbe.InotI} \textbf{(a)}.

\textbf{(b)} Clearly,
\begin{equation}
\sum I=\sum_{i\in I}i=\sum_{x\in I}x \label{pf.lem.sol.Ialbe.InotI.b.triv}%
\end{equation}
(here, we have renamed the summation index $i$ as $x$).

We know that $\left(  a_{1},a_{2},\ldots,a_{k}\right)  $ is a list of all
elements of $I$ (with no repetitions). Thus, Lemma \ref{lem.sol.Ialbe.inclist}
\textbf{(a)} (applied to $I$, $k$ and $\left(  a_{1},a_{2},\ldots
,a_{k}\right)  $ instead of $S$, $s$ and $\left(  c_{1},c_{2},\ldots
,c_{s}\right)  $) shows that the map $\left[  k\right]  \rightarrow
I,\ h\mapsto a_{h}$ is well-defined and a bijection.

We know that $\left(  b_{1},b_{2},\ldots,b_{n-k}\right)  $ is a list of all
elements of $\left[  n\right]  \setminus I$ (with no repetitions). Thus, Lemma
\ref{lem.sol.Ialbe.inclist} \textbf{(a)} (applied to $\left[  n\right]
\setminus I$, $n-k$ and $\left(  b_{1},b_{2},\ldots,b_{n-k}\right)  $ instead
of $S$, $s$ and $\left(  c_{1},c_{2},\ldots,c_{s}\right)  $) shows that the
map $\left[  n-k\right]  \rightarrow\left[  n\right]  \setminus I,\ h\mapsto
b_{h}$ is well-defined and a bijection.

Now,%
\begin{align}
&  \underbrace{\sum_{\left(  i,j\right)  \in\left[  k\right]  \times\left[
n-k\right]  }}_{=\sum_{i\in\left[  k\right]  }\sum_{j\in\left[  n-k\right]  }%
}\left[  a_{i}>b_{j}\right] \nonumber\\
&  =\sum_{i\in\left[  k\right]  }\underbrace{\sum_{j\in\left[  n-k\right]
}\left[  a_{i}>b_{j}\right]  }_{\substack{=\sum_{h\in\left[  n-k\right]
}\left[  a_{i}>b_{h}\right]  \\\text{(here, we have renamed the}%
\\\text{summation index }j\text{ as }h\text{)}}}=\sum_{i\in\left[  k\right]
}\underbrace{\sum_{h\in\left[  n-k\right]  }\left[  a_{i}>b_{h}\right]
}_{\substack{=\sum_{y\in\left[  n\right]  \setminus I}\left[  a_{i}>y\right]
\\\text{(here, we have substituted }y\text{ for }b_{h}\text{ in the
sum,}\\\text{since the map }\left[  n-k\right]  \rightarrow\left[  n\right]
\setminus I,\ h\mapsto b_{h}\text{ is a bijection)}}}\nonumber\\
&  =\sum_{i\in\left[  k\right]  }\sum_{y\in\left[  n\right]  \setminus
I}\left[  a_{i}>y\right] \nonumber\\
&  =\sum_{h\in\left[  k\right]  }\sum_{y\in\left[  n\right]  \setminus
I}\left[  a_{h}>y\right]  \ \ \ \ \ \ \ \ \ \ \left(
\begin{array}
[c]{c}%
\text{here, we have renamed the}\\
\text{summation index }i\text{ as }h\text{ in the outer sum}%
\end{array}
\right) \nonumber\\
&  =\sum_{x\in I}\sum_{y\in\left[  n\right]  \setminus I}\left[  x>y\right]
\label{pf.lem.sol.Ialbe.InotI.b.1}%
\end{align}
(here, we have substituted $x$ for $a_{h}$ in the outer sum, since the map
$\left[  k\right]  \rightarrow I,\ h\mapsto a_{h}$ is a bijection). But every
$x\in I$ satisfies%
\begin{equation}
\sum_{y\in\left[  n\right]  \setminus I}\left[  x>y\right]  =x-1-\sum_{y\in
I}\left[  x>y\right]  \label{pf.lem.sol.Ialbe.InotI.b.sum1}%
\end{equation}
\footnote{\textit{Proof of (\ref{pf.lem.sol.Ialbe.InotI.b.sum1}):} Let $x\in
I$. Then, $x\in I\subseteq\left[  n\right]  =\left\{  1,2,\ldots,n\right\}  $,
so that $1\leq x\leq n$.
\par
Every $y\in\left[  n\right]  $ satisfies either $y\in I$ or $y\notin I$ (but
not both). Hence,%
\begin{align*}
\sum_{y\in\left[  n\right]  }\left[  x>y\right]   &  =\underbrace{\sum
_{\substack{y\in\left[  n\right]  ;\\y\in I}}}_{\substack{=\sum_{y\in
I}\\\text{(since }I\subseteq\left[  n\right]  \text{)}}}\left[  x>y\right]
+\underbrace{\sum_{\substack{y\in\left[  n\right]  ;\\y\notin I}}}%
_{=\sum_{y\in\left[  n\right]  \setminus I}}\left[  x>y\right] \\
&  =\sum_{y\in I}\left[  x>y\right]  +\sum_{y\in\left[  n\right]  \setminus
I}\left[  x>y\right]  .
\end{align*}
Comparing this with%
\begin{align*}
\underbrace{\sum_{y\in\left[  n\right]  }}_{\substack{=\sum_{y\in\left\{
1,2,\ldots,n\right\}  }\\\text{(since }\left[  n\right]  =\left\{
1,2,\ldots,n\right\}  \text{)}}}\left[  x>y\right]   &  =\underbrace{\sum
_{y\in\left\{  1,2,\ldots,n\right\}  }}_{=\sum_{y=1}^{n}}\left[  x>y\right]
=\sum_{y=1}^{n}\left[  x>y\right] \\
&  =\sum_{y=1}^{x-1}\underbrace{\left[  x>y\right]  }%
_{\substack{=1\\\text{(since }x>y\\\text{(since }y\leq x-1<x\text{))}}%
}+\sum_{y=x}^{n}\underbrace{\left[  x>y\right]  }_{\substack{=0\\\text{(since
we don't have }x>y\\\text{(since }y\geq x\text{))}}%
}\ \ \ \ \ \ \ \ \ \ \left(  \text{since }1\leq x\leq n\right) \\
&  =\sum_{y=1}^{x-1}1+\underbrace{\sum_{y=x}^{n}0}_{=0}=\sum_{y=1}%
^{x-1}1=\left(  x-1\right)  1=x-1,
\end{align*}
we obtain%
\[
\sum_{y\in I}\left[  x>y\right]  +\sum_{y\in\left[  n\right]  \setminus
I}\left[  x>y\right]  =x-1.
\]
Subtracting $\sum_{y\in I}\left[  x>y\right]  $ from both sides of this
equality, we obtain $\sum_{y\in\left[  n\right]  \setminus I}\left[
x>y\right]  =x-1-\sum_{y\in I}\left[  x>y\right]  $. This proves
(\ref{pf.lem.sol.Ialbe.InotI.b.sum1}).}. Thus,
(\ref{pf.lem.sol.Ialbe.InotI.b.1}) becomes%
\begin{align*}
&  \sum_{\left(  i,j\right)  \in\left[  k\right]  \times\left[  n-k\right]
}\left[  a_{i}>b_{j}\right] \\
&  =\sum_{x\in I}\underbrace{\sum_{y\in\left[  n\right]  \setminus I}\left[
x>y\right]  }_{\substack{=x-1-\sum_{y\in I}\left[  x>y\right]  \\\text{(by
(\ref{pf.lem.sol.Ialbe.InotI.b.sum1}))}}}\\
&  =\sum_{x\in I}\left(  x-1-\sum_{y\in I}\left[  x>y\right]  \right)
=\underbrace{\sum_{x\in I}x}_{\substack{=\sum I\\\text{(by
(\ref{pf.lem.sol.Ialbe.InotI.b.triv}))}}}-\underbrace{\sum_{x\in I}%
1}_{=\left\vert I\right\vert \cdot1=\left\vert I\right\vert =k}%
-\underbrace{\sum_{x\in I}\sum_{y\in I}\left[  x>y\right]  }%
_{\substack{=0+1+\cdots+\left(  k-1\right)  \\\text{(by Lemma
\ref{lem.sol.Ialbe.II})}}}\\
&  =\sum I-k-\left(  0+1+\cdots+\left(  k-1\right)  \right)  =\sum
I-\underbrace{\left(  \left(  0+1+\cdots+\left(  k-1\right)  \right)
+k\right)  }_{=0+1+\cdots+k=1+2+\cdots+k}\\
&  =\sum I-\left(  1+2+\cdots+k\right)  .
\end{align*}
This proves Lemma \ref{lem.sol.Ialbe.InotI} \textbf{(b)}.
\end{proof}
\end{verlong}

Let us now come to the actual solution of Exercise \ref{exe.Ialbe}:

\begin{proof}
[Solution to Exercise \ref{exe.Ialbe}.]We have $I\subseteq\left\{
1,2,\ldots,n\right\}  =\left[  n\right]  $ (since $\left[  n\right]  =\left\{
1,2,\ldots,n\right\}  $). Hence, $\left\vert \left[  n\right]  \setminus
I\right\vert =\left\vert \underbrace{\left[  n\right]  }_{=\left\{
1,2,\ldots,n\right\}  }\right\vert -\underbrace{\left\vert I\right\vert }%
_{=k}=\underbrace{\left\vert \left\{  1,2,\ldots,n\right\}  \right\vert }%
_{=n}-k=n-k$. Thus, $n-k=\left\vert \left[  n\right]  \setminus I\right\vert
\in\mathbb{N}$ (since the cardinality of any finite set is $\in\mathbb{N}$).
Hence, $n-k\geq0$, so that $k\leq n$. Also, $k=\left\vert I\right\vert
\in\mathbb{N}$.

We know that $\left(  a_{1},a_{2},\ldots,a_{k}\right)  $ is a list of all
elements of $I$ (with no repetitions). Thus, Lemma \ref{lem.sol.Ialbe.inclist}
\textbf{(b)} (applied to $I$, $k$, $\left(  a_{1},a_{2},\ldots,a_{k}\right)  $
and $\alpha$ instead of $S$, $s$, $\left(  c_{1},c_{2},\ldots,c_{s}\right)  $
and $\pi$) shows that $\left(  a_{\alpha\left(  1\right)  },a_{\alpha\left(
2\right)  },\ldots,a_{\alpha\left(  k\right)  }\right)  $ is a list of all
elements of $I$ (with no repetitions).

We know that $\left(  b_{1},b_{2},\ldots,b_{n-k}\right)  $ is a list of all
elements of $\left\{  1,2,\ldots,n\right\}  \setminus I$ (with no
repetitions). In other words, $\left(  b_{1},b_{2},\ldots,b_{n-k}\right)  $ is
a list of all elements of $\left[  n\right]  \setminus I$ (with no
repetitions) (since $\left[  n\right]  =\left\{  1,2,\ldots,n\right\}  $).
Thus, Lemma \ref{lem.sol.Ialbe.inclist} \textbf{(b)} (applied to $\left[
n\right]  \setminus I$, $n-k$, $\left(  b_{1},b_{2},\ldots,b_{n-k}\right)  $
and $\beta$ instead of $S$, $s$, $\left(  c_{1},c_{2},\ldots,c_{s}\right)  $
and $\pi$) shows that $\left(  b_{\beta\left(  1\right)  },b_{\beta\left(
2\right)  },\ldots,b_{\beta\left(  n-k\right)  }\right)  $ is a list of all
elements of $\left[  n\right]  \setminus I$ (with no repetitions).

\textbf{(a)} Lemma \ref{lem.sol.Ialbe.InotI} \textbf{(a)} (applied to $\left(
a_{\alpha\left(  1\right)  },a_{\alpha\left(  2\right)  },\ldots
,a_{\alpha\left(  k\right)  }\right)  $ and $\left(  b_{\beta\left(  1\right)
},b_{\beta\left(  2\right)  },\ldots,b_{\beta\left(  n-k\right)  }\right)  $
instead of $\left(  a_{1},a_{2},\ldots,a_{k}\right)  $ and $\left(
b_{1},b_{2},\ldots,b_{n-k}\right)  $) shows that there exists a unique
$\sigma\in S_{n}$ satisfying%
\[
\left(  \sigma\left(  1\right)  ,\sigma\left(  2\right)  ,\ldots,\sigma\left(
n\right)  \right)  =\left(  a_{\alpha\left(  1\right)  },a_{\alpha\left(
2\right)  },\ldots,a_{\alpha\left(  k\right)  },b_{\beta\left(  1\right)
},b_{\beta\left(  2\right)  },\ldots,b_{\beta\left(  n-k\right)  }\right)  .
\]
This solves Exercise \ref{exe.Ialbe} \textbf{(a)}.

\textbf{(b)} Let $\mathbf{a}$ be the $k$-tuple $\left(  a_{\alpha\left(
1\right)  },a_{\alpha\left(  2\right)  },\ldots,a_{\alpha\left(  k\right)
}\right)  $ of integers. Recall that $\left(  a_{1},a_{2},\ldots,a_{k}\right)
$ is the list of all elements of $I$ in increasing order (with no
repetitions). Moreover, $k=\left\vert I\right\vert $. Hence, Lemma
\ref{lem.sol.Ialbe.Inv=Inv} \textbf{(b)} (applied to $I$, $k$, $\alpha$ and
$\left(  a_{1},a_{2},\ldots,a_{k}\right)  $ instead of $P$, $m$, $\sigma$ and
$\left(  p_{1},p_{2},\ldots,p_{m}\right)  $) yields%
\begin{equation}
\ell\left(  \alpha\right)  =\ell\left(  \underbrace{\left(  a_{\alpha\left(
1\right)  },a_{\alpha\left(  2\right)  },\ldots,a_{\alpha\left(  k\right)
}\right)  }_{=\mathbf{a}}\right)  =\ell\left(  \mathbf{a}\right)  .
\label{sol.Ialbe.b.la}%
\end{equation}

Let $\mathbf{b}$ be the $\left(  n-k\right)  $-tuple $\left(  b_{\beta\left(
1\right)  },b_{\beta\left(  2\right)  },\ldots,b_{\beta\left(  n-k\right)
}\right)  $ of integers. Recall that \newline$\left(  b_{1},b_{2}%
,\ldots,b_{n-k}\right)  $ is the list of all elements of $\left\{
1,2,\ldots,n\right\}  \setminus I$ in increasing order (with no repetitions).
In other words, $\left(  b_{1},b_{2},\ldots,b_{n-k}\right)  $ is the list of
all elements of $\left[  n\right]  \setminus I$ in increasing order (with no
repetitions) (since $\left[  n\right]  =\left\{  1,2,\ldots,n\right\}  $).
Moreover, $n-k=\left\vert \left[  n\right]  \setminus I\right\vert $ (since
$\left\vert \left[  n\right]  \setminus I\right\vert =n-k$). Hence, Lemma
\ref{lem.sol.Ialbe.Inv=Inv} \textbf{(b)} (applied to $\left[  n\right]
\setminus I$, $n-k$, $\beta$ and $\left(  b_{1},b_{2},\ldots,b_{n-k}\right)  $
instead of $P$, $m$, $\sigma$ and $\left(  p_{1},p_{2},\ldots,p_{m}\right)  $)
yields%
\begin{equation}
\ell\left(  \beta\right)  =\ell\left(  \underbrace{\left(  b_{\beta\left(
1\right)  },b_{\beta\left(  2\right)  },\ldots,b_{\beta\left(  n-k\right)
}\right)  }_{=\mathbf{b}}\right)  =\ell\left(  \mathbf{b}\right)  .
\label{sol.Ialbe.b.lb}%
\end{equation}

Let $\mathbf{c}$ be the $\left(  k+\left(  n-k\right)  \right)  $-tuple
$\left(  a_{\alpha\left(  1\right)  },a_{\alpha\left(  2\right)  }%
,\ldots,a_{\alpha\left(  k\right)  },b_{\beta\left(  1\right)  }%
,b_{\beta\left(  2\right)  },\ldots,b_{\beta\left(  n-k\right)  }\right)  $ of integers.

The permutation $\sigma_{I,\alpha,\beta}$ is the unique $\sigma\in S_{n}$
satisfying%
\begin{equation}
\left(  \sigma\left(  1\right)  ,\sigma\left(  2\right)  ,\ldots,\sigma\left(
n\right)  \right)  =\left(  a_{\alpha\left(  1\right)  },a_{\alpha\left(
2\right)  },\ldots,a_{\alpha\left(  k\right)  },b_{\beta\left(  1\right)
},b_{\beta\left(  2\right)  },\ldots,b_{\beta\left(  n-k\right)  }\right)  .
\label{sol.Ialbe.b.1}%
\end{equation}
Thus, $\sigma_{I,\alpha,\beta}$ is a $\sigma\in S_{n}$ satisfying
(\ref{sol.Ialbe.b.1}). In other words, we have%
\begin{align}
&  \left(  \sigma_{I,\alpha,\beta}\left(  1\right)  ,\sigma_{I,\alpha,\beta
}\left(  2\right)  ,\ldots,\sigma_{I,\alpha,\beta}\left(  n\right)  \right)
\nonumber\\
&  =\left(  a_{\alpha\left(  1\right)  },a_{\alpha\left(  2\right)  }%
,\ldots,a_{\alpha\left(  k\right)  },b_{\beta\left(  1\right)  }%
,b_{\beta\left(  2\right)  },\ldots,b_{\beta\left(  n-k\right)  }\right)
=\mathbf{c} \label{sol.Ialbe.b.2}%
\end{align}
(since $\mathbf{c}=\left(  a_{\alpha\left(  1\right)  },a_{\alpha\left(
2\right)  },\ldots,a_{\alpha\left(  k\right)  },b_{\beta\left(  1\right)
},b_{\beta\left(  2\right)  },\ldots,b_{\beta\left(  n-k\right)  }\right)  $).

Now, $\left\{  1,2,\ldots,n\right\}  $ is a finite set of integers and
satisfies $n=\left\vert \left\{  1,2,\ldots,n\right\}  \right\vert $.
Moreover, $\left(  1,2,\ldots,n\right)  $ is the list of all elements of
$\left\{  1,2,\ldots,n\right\}  $ in increasing order (with no repetitions).
Hence, Lemma \ref{lem.sol.Ialbe.Inv=Inv} \textbf{(b)} (applied to $\left\{
1,2,\ldots,n\right\}  $, $n$, $\sigma_{I,\alpha,\beta}$ and $\left(
1,2,\ldots,n\right)  $) yields%
\begin{align}
\ell\left(  \sigma_{I,\alpha,\beta}\right)   &  =\ell\left(
\underbrace{\left(  \sigma_{I,\alpha,\beta}\left(  1\right)  ,\sigma
_{I,\alpha,\beta}\left(  2\right)  ,\ldots,\sigma_{I,\alpha,\beta}\left(
n\right)  \right)  }_{\substack{=\mathbf{c}\\\text{(by (\ref{sol.Ialbe.b.2}%
))}}}\right) \nonumber\\
&  =\ell\left(  \mathbf{c}\right)  =\ell\left(  \mathbf{a}\right)
+\ell\left(  \mathbf{b}\right)  +\sum_{\left(  i,j\right)  \in\left[
k\right]  \times\left[  n-k\right]  }\left[  a_{\alpha\left(  i\right)
}>b_{\beta\left(  j\right)  }\right]  \label{sol.Ialbe.b.9}%
\end{align}
(by Lemma \ref{lem.sol.Ialbe.ab} (applied to $k$, $n-k$, $\left(
a_{\alpha\left(  1\right)  },a_{\alpha\left(  2\right)  },\ldots
,a_{\alpha\left(  k\right)  }\right)  $ and $\left(  b_{\beta\left(  1\right)
},b_{\beta\left(  2\right)  },\ldots,b_{\beta\left(  n-k\right)  }\right)  $
instead of $n$, $m$, $\left(  a_{1},a_{2},\ldots,a_{n}\right)  $ and $\left(
b_{1},b_{2},\ldots,b_{m}\right)  $)).

Now, recall that $\left(  a_{\alpha\left(  1\right)  },a_{\alpha\left(
2\right)  },\ldots,a_{\alpha\left(  k\right)  }\right)  $ is a list of all
elements of $I$ (with no repetitions). Also, recall that $\left(
b_{\beta\left(  1\right)  },b_{\beta\left(  2\right)  },\ldots,b_{\beta\left(
n-k\right)  }\right)  $ is a list of all elements of $\left[  n\right]
\setminus I$ (with no repetitions). Lemma \ref{lem.sol.Ialbe.InotI}
\textbf{(b)} (applied to $\left(  a_{\alpha\left(  1\right)  },a_{\alpha
\left(  2\right)  },\ldots,a_{\alpha\left(  k\right)  }\right)  $ and $\left(
b_{\beta\left(  1\right)  },b_{\beta\left(  2\right)  },\ldots,b_{\beta\left(
n-k\right)  }\right)  $ instead of $\left(  a_{1},a_{2},\ldots,a_{k}\right)  $
and $\left(  b_{1},b_{2},\ldots,b_{n-k}\right)  $) thus shows that%
\[
\sum_{\left(  i,j\right)  \in\left[  k\right]  \times\left[  n-k\right]
}\left[  a_{\alpha\left(  i\right)  }>b_{\beta\left(  j\right)  }\right]
=\sum I-\left(  1+2+\cdots+k\right)  .
\]
Hence, (\ref{sol.Ialbe.b.9}) becomes%
\begin{align*}
\ell\left(  \sigma_{I,\alpha,\beta}\right)   &  =\underbrace{\ell\left(
\mathbf{a}\right)  }_{\substack{=\ell\left(  \alpha\right)  \\\text{(by
(\ref{sol.Ialbe.b.la}))}}}+\underbrace{\ell\left(  \mathbf{b}\right)
}_{\substack{=\ell\left(  \beta\right)  \\\text{(by (\ref{sol.Ialbe.b.lb}))}%
}}+\underbrace{\sum_{\left(  i,j\right)  \in\left[  k\right]  \times\left[
n-k\right]  }\left[  a_{\alpha\left(  i\right)  }>b_{\beta\left(  j\right)
}\right]  }_{=\sum I-\left(  1+2+\cdots+k\right)  }\\
&  =\ell\left(  \alpha\right)  +\ell\left(  \beta\right)  +\sum I-\left(
1+2+\cdots+k\right)  .
\end{align*}

It thus remains to prove that $\left(  -1\right)  ^{\sigma_{I,\alpha,\beta}%
}=\left(  -1\right)  ^{\alpha}\cdot\left(  -1\right)  ^{\beta}\cdot\left(
-1\right)  ^{\sum I-\left(  1+2+\cdots+k\right)  }$.

But the definition of $\left(  -1\right)  ^{\sigma_{I,\alpha,\beta}}$ yields
\begin{align*}
\left(  -1\right)  ^{\sigma_{I,\alpha,\beta}}  &  =\left(  -1\right)
^{\ell\left(  \sigma_{I,\alpha,\beta}\right)  }=\left(  -1\right)
^{\ell\left(  \alpha\right)  +\ell\left(  \beta\right)  +\sum I-\left(
1+2+\cdots+k\right)  }\\
&  \ \ \ \ \ \ \ \ \ \ \left(  \text{since }\ell\left(  \sigma_{I,\alpha
,\beta}\right)  =\ell\left(  \alpha\right)  +\ell\left(  \beta\right)  +\sum
I-\left(  1+2+\cdots+k\right)  \right) \\
&  =\left(  -1\right)  ^{\ell\left(  \alpha\right)  }\cdot\left(  -1\right)
^{\ell\left(  \beta\right)  }\cdot\left(  -1\right)  ^{\sum I-\left(
1+2+\cdots+k\right)  }.
\end{align*}
Comparing this with%
\begin{align*}
&  \underbrace{\left(  -1\right)  ^{\alpha}}_{\substack{=\left(  -1\right)
^{\ell\left(  \alpha\right)  }\\\text{(by the definition of }\left(
-1\right)  ^{\alpha}\text{)}}}\cdot\underbrace{\left(  -1\right)  ^{\beta}%
}_{\substack{=\left(  -1\right)  ^{\ell\left(  \beta\right)  }\\\text{(by the
definition of }\left(  -1\right)  ^{\beta}\text{)}}}\cdot\left(  -1\right)
^{\sum I-\left(  1+2+\cdots+k\right)  }\\
&  =\left(  -1\right)  ^{\ell\left(  \alpha\right)  }\cdot\left(  -1\right)
^{\ell\left(  \beta\right)  }\cdot\left(  -1\right)  ^{\sum I-\left(
1+2+\cdots+k\right)  },
\end{align*}
this yields $\left(  -1\right)  ^{\sigma_{I,\alpha,\beta}}=\left(  -1\right)
^{\alpha}\cdot\left(  -1\right)  ^{\beta}\cdot\left(  -1\right)  ^{\sum
I-\left(  1+2+\cdots+k\right)  }$. This completes the solution of Exercise
\ref{exe.Ialbe} \textbf{(b)}.

\begin{vershort}
\textbf{(c)} For every $\left(  \alpha,\beta\right)  \in S_{k}\times S_{n-k}$,
the permutation $\sigma_{I,\alpha,\beta}\in S_{n}$ satisfies%
\begin{align}
&  \left(  \sigma_{I,\alpha,\beta}\left(  1\right)  ,\sigma_{I,\alpha,\beta
}\left(  2\right)  ,\ldots,\sigma_{I,\alpha,\beta}\left(  n\right)  \right)
\nonumber\\
&  =\left(  a_{\alpha\left(  1\right)  },a_{\alpha\left(  2\right)  }%
,\ldots,a_{\alpha\left(  k\right)  },b_{\beta\left(  1\right)  }%
,b_{\beta\left(  2\right)  },\ldots,b_{\beta\left(  n-k\right)  }\right)
\label{sol.Ialbe.c.short.siIAB}%
\end{align}
\footnote{\textit{Proof of (\ref{sol.Ialbe.c.short.siIAB}):} Let $\left(
\alpha,\beta\right)  \in S_{k}\times S_{n-k}$. The permutation $\sigma
_{I,\alpha,\beta}$ is the unique $\sigma\in S_{n}$ satisfying%
\begin{equation}
\left(  \sigma\left(  1\right)  ,\sigma\left(  2\right)  ,\ldots,\sigma\left(
n\right)  \right)  =\left(  a_{\alpha\left(  1\right)  },a_{\alpha\left(
2\right)  },\ldots,a_{\alpha\left(  k\right)  },b_{\beta\left(  1\right)
},b_{\beta\left(  2\right)  },\ldots,b_{\beta\left(  n-k\right)  }\right)  .
\label{sol.Ialbe.c.short.siIAB.pf.1}%
\end{equation}
Thus, $\sigma_{I,\alpha,\beta}$ is a $\sigma\in S_{n}$ satisfying
(\ref{sol.Ialbe.c.short.siIAB.pf.1}). In other words, we have%
\[
\left(  \sigma_{I,\alpha,\beta}\left(  1\right)  ,\sigma_{I,\alpha,\beta
}\left(  2\right)  ,\ldots,\sigma_{I,\alpha,\beta}\left(  n\right)  \right)
=\left(  a_{\alpha\left(  1\right)  },a_{\alpha\left(  2\right)  }%
,\ldots,a_{\alpha\left(  k\right)  },b_{\beta\left(  1\right)  }%
,b_{\beta\left(  2\right)  },\ldots,b_{\beta\left(  n-k\right)  }\right)  .
\]
Qed.}.

For every $\left(  \alpha,\beta\right)  \in S_{k}\times S_{n-k}$, the element
$\sigma_{I,\alpha,\beta}$ is a well-defined element of $\left\{  \tau\in
S_{n}\ \mid\ \tau\left(  \left\{  1,2,\ldots,k\right\}  \right)  =I\right\}
$\ \ \ \ \footnote{\textit{Proof.} Let $\left(  \alpha,\beta\right)  \in
S_{k}\times S_{n-k}$. We must show that the element $\sigma_{I,\alpha,\beta}$
is a well-defined element of $\left\{  \tau\in S_{n}\ \mid\ \tau\left(
\left\{  1,2,\ldots,k\right\}  \right)  =I\right\}  $.
\par
Now,
\begin{align*}
&  \left(  \sigma_{I,\alpha,\beta}\left(  1\right)  ,\sigma_{I,\alpha,\beta
}\left(  2\right)  ,\ldots,\sigma_{I,\alpha,\beta}\left(  k\right)  \right) \\
&  =\left(  \text{the list of the first }k\text{ entries of the list
}\underbrace{\left(  \sigma_{I,\alpha,\beta}\left(  1\right)  ,\sigma
_{I,\alpha,\beta}\left(  2\right)  ,\ldots,\sigma_{I,\alpha,\beta}\left(
n\right)  \right)  }_{\substack{=\left(  a_{\alpha\left(  1\right)
},a_{\alpha\left(  2\right)  },\ldots,a_{\alpha\left(  k\right)  }%
,b_{\beta\left(  1\right)  },b_{\beta\left(  2\right)  },\ldots,b_{\beta
\left(  n-k\right)  }\right)  \\\text{(by (\ref{sol.Ialbe.c.short.siIAB}))}%
}}\right) \\
&  =\left(  \text{the list of the first }k\text{ entries of the list }\left(
a_{\alpha\left(  1\right)  },a_{\alpha\left(  2\right)  },\ldots
,a_{\alpha\left(  k\right)  },b_{\beta\left(  1\right)  },b_{\beta\left(
2\right)  },\ldots,b_{\beta\left(  n-k\right)  }\right)  \right) \\
&  =\left(  a_{\alpha\left(  1\right)  },a_{\alpha\left(  2\right)  }%
,\ldots,a_{\alpha\left(  k\right)  }\right)  .
\end{align*}
Hence,%
\[
\left\{  \sigma_{I,\alpha,\beta}\left(  1\right)  ,\sigma_{I,\alpha,\beta
}\left(  2\right)  ,\ldots,\sigma_{I,\alpha,\beta}\left(  k\right)  \right\}
=\left\{  a_{\alpha\left(  1\right)  },a_{\alpha\left(  2\right)  }%
,\ldots,a_{\alpha\left(  k\right)  }\right\}  =I
\]
(since $\left(  a_{\alpha\left(  1\right)  },a_{\alpha\left(  2\right)
},\ldots,a_{\alpha\left(  k\right)  }\right)  $ is a list of all elements of
$I$). Now,%
\[
\sigma_{I,\alpha,\beta}\left(  \left\{  1,2,\ldots,k\right\}  \right)
=\left\{  \sigma_{I,\alpha,\beta}\left(  1\right)  ,\sigma_{I,\alpha,\beta
}\left(  2\right)  ,\ldots,\sigma_{I,\alpha,\beta}\left(  k\right)  \right\}
=I.
\]
\par
So we know that $\sigma_{I,\alpha,\beta}$ is an element of $S_{n}$ and
satisfies $\sigma_{I,\alpha,\beta}\left(  \left\{  1,2,\ldots,k\right\}
\right)  =I$. In other words,%
\[
\sigma_{I,\alpha,\beta}\in\left\{  \tau\in S_{n}\ \mid\ \tau\left(  \left\{
1,2,\ldots,k\right\}  \right)  =I\right\}  .
\]
We thus have proven that $\sigma_{I,\alpha,\beta}$ is a well-defined element
of $\left\{  \tau\in S_{n}\ \mid\ \tau\left(  \left\{  1,2,\ldots,k\right\}
\right)  =I\right\}  $. Qed.}. Hence, the map%
\begin{align*}
S_{k}\times S_{n-k}  &  \rightarrow\left\{  \tau\in S_{n}\ \mid\ \tau\left(
\left\{  1,2,\ldots,k\right\}  \right)  =I\right\}  ,\\
\left(  \alpha,\beta\right)   &  \mapsto\sigma_{I,\alpha,\beta}%
\end{align*}
is well-defined. Denote this map by $\mu$.
\end{vershort}

\begin{verlong}
\textbf{(c)} For every $\left(  \alpha,\beta\right)  \in S_{k}\times S_{n-k}$,
the element $\sigma_{I,\alpha,\beta}$ is a well-defined element of $\left\{
\tau\in S_{n}\ \mid\ \tau\left(  \left\{  1,2,\ldots,k\right\}  \right)
=I\right\}  $\ \ \ \ \footnote{\textit{Proof.} Let $\left(  \alpha
,\beta\right)  \in S_{k}\times S_{n-k}$. We must show that the element
$\sigma_{I,\alpha,\beta}$ is a well-defined element of $\left\{  \tau\in
S_{n}\ \mid\ \tau\left(  \left\{  1,2,\ldots,k\right\}  \right)  =I\right\}
$.
\par
We have $\left(  \alpha,\beta\right)  \in S_{k}\times S_{n-k}$. In other
words, $\alpha\in S_{k}$ and $\beta\in S_{n-k}$. Hence, $\sigma_{I,\alpha
,\beta}$ is a well-defined element of $S_{n}$.
\par
The permutation $\sigma_{I,\alpha,\beta}$ is the unique $\sigma\in S_{n}$
satisfying%
\begin{equation}
\left(  \sigma\left(  1\right)  ,\sigma\left(  2\right)  ,\ldots,\sigma\left(
n\right)  \right)  =\left(  a_{\alpha\left(  1\right)  },a_{\alpha\left(
2\right)  },\ldots,a_{\alpha\left(  k\right)  },b_{\beta\left(  1\right)
},b_{\beta\left(  2\right)  },\ldots,b_{\beta\left(  n-k\right)  }\right)  .
\label{sol.Ialbe.c.wd.fn1.1}%
\end{equation}
Thus, $\sigma_{I,\alpha,\beta}$ is a $\sigma\in S_{n}$ satisfying
(\ref{sol.Ialbe.c.wd.fn1.1}). In other words, we have%
\begin{align*}
&  \left(  \sigma_{I,\alpha,\beta}\left(  1\right)  ,\sigma_{I,\alpha,\beta
}\left(  2\right)  ,\ldots,\sigma_{I,\alpha,\beta}\left(  n\right)  \right) \\
&  =\left(  a_{\alpha\left(  1\right)  },a_{\alpha\left(  2\right)  }%
,\ldots,a_{\alpha\left(  k\right)  },b_{\beta\left(  1\right)  }%
,b_{\beta\left(  2\right)  },\ldots,b_{\beta\left(  n-k\right)  }\right)  .
\end{align*}
Now,
\begin{align*}
&  \left(  \sigma_{I,\alpha,\beta}\left(  1\right)  ,\sigma_{I,\alpha,\beta
}\left(  2\right)  ,\ldots,\sigma_{I,\alpha,\beta}\left(  k\right)  \right) \\
&  =\left(  \text{the list of the first }k\text{ entries of the list
}\underbrace{\left(  \sigma_{I,\alpha,\beta}\left(  1\right)  ,\sigma
_{I,\alpha,\beta}\left(  2\right)  ,\ldots,\sigma_{I,\alpha,\beta}\left(
n\right)  \right)  }_{=\left(  a_{\alpha\left(  1\right)  },a_{\alpha\left(
2\right)  },\ldots,a_{\alpha\left(  k\right)  },b_{\beta\left(  1\right)
},b_{\beta\left(  2\right)  },\ldots,b_{\beta\left(  n-k\right)  }\right)
}\right) \\
&  =\left(  \text{the list of the first }k\text{ entries of the list }\left(
a_{\alpha\left(  1\right)  },a_{\alpha\left(  2\right)  },\ldots
,a_{\alpha\left(  k\right)  },b_{\beta\left(  1\right)  },b_{\beta\left(
2\right)  },\ldots,b_{\beta\left(  n-k\right)  }\right)  \right) \\
&  =\left(  a_{\alpha\left(  1\right)  },a_{\alpha\left(  2\right)  }%
,\ldots,a_{\alpha\left(  k\right)  }\right)  .
\end{align*}
Hence,%
\[
\left\{  \sigma_{I,\alpha,\beta}\left(  1\right)  ,\sigma_{I,\alpha,\beta
}\left(  2\right)  ,\ldots,\sigma_{I,\alpha,\beta}\left(  k\right)  \right\}
=\left\{  a_{\alpha\left(  1\right)  },a_{\alpha\left(  2\right)  }%
,\ldots,a_{\alpha\left(  k\right)  }\right\}  =I
\]
(since $\left(  a_{\alpha\left(  1\right)  },a_{\alpha\left(  2\right)
},\ldots,a_{\alpha\left(  k\right)  }\right)  $ is a list of all elements of
$I$ (with no repetitions)). Now,%
\[
\sigma_{I,\alpha,\beta}\left(  \left\{  1,2,\ldots,k\right\}  \right)
=\left\{  \sigma_{I,\alpha,\beta}\left(  1\right)  ,\sigma_{I,\alpha,\beta
}\left(  2\right)  ,\ldots,\sigma_{I,\alpha,\beta}\left(  k\right)  \right\}
=I.
\]
\par
So we know that $\sigma_{I,\alpha,\beta}$ is an element of $S_{n}$ and
satisfies $\sigma_{I,\alpha,\beta}\left(  \left\{  1,2,\ldots,k\right\}
\right)  =I$. In other words, $\sigma_{I,\alpha,\beta}$ is an element $\tau\in
S_{n}$ satisfying $\tau\left(  \left\{  1,2,\ldots,k\right\}  \right)  =I$. In
other words,%
\[
\sigma_{I,\alpha,\beta}\in\left\{  \tau\in S_{n}\ \mid\ \tau\left(  \left\{
1,2,\ldots,k\right\}  \right)  =I\right\}  .
\]
We thus have proven that $\sigma_{I,\alpha,\beta}$ is a well-defined element
of $\left\{  \tau\in S_{n}\ \mid\ \tau\left(  \left\{  1,2,\ldots,k\right\}
\right)  =I\right\}  $. Qed.}. Hence, the map%
\begin{align*}
S_{k}\times S_{n-k}  &  \rightarrow\left\{  \tau\in S_{n}\ \mid\ \tau\left(
\left\{  1,2,\ldots,k\right\}  \right)  =I\right\}  ,\\
\left(  \alpha,\beta\right)   &  \mapsto\sigma_{I,\alpha,\beta}%
\end{align*}
is well-defined. Denote this map by $\mu$.
\end{verlong}

\begin{vershort}
Now, it is easy to see that the map $\mu$ is
injective\footnote{\textit{Proof.} Let $\left(  \alpha,\beta\right)  $ and
$\left(  \alpha^{\prime},\beta^{\prime}\right)  $ be two elements of
$S_{k}\times S_{n-k}$ satisfying $\mu\left(  \alpha,\beta\right)  =\mu\left(
\alpha^{\prime},\beta^{\prime}\right)  $. We will show that $\left(
\alpha,\beta\right)  =\left(  \alpha^{\prime},\beta^{\prime}\right)  $.
\par
The definition of $\mu\left(  \alpha,\beta\right)  $ yields $\mu\left(
\alpha,\beta\right)  =\sigma_{I,\alpha,\beta}$. The definition of $\mu\left(
\alpha^{\prime},\beta^{\prime}\right)  $ yields $\mu\left(  \alpha^{\prime
},\beta^{\prime}\right)  =\sigma_{I,\alpha^{\prime},\beta^{\prime}}$. Hence,
$\sigma_{I,\alpha^{\prime},\beta^{\prime}}=\mu\left(  \alpha^{\prime}%
,\beta^{\prime}\right)  =\mu\left(  \alpha,\beta\right)  =\sigma
_{I,\alpha,\beta}$.
\par
Now, (\ref{sol.Ialbe.c.short.siIAB}) (applied to $\left(  \alpha^{\prime
},\beta^{\prime}\right)  $ instead of $\left(  \alpha,\beta\right)  $) yields%
\begin{align*}
&  \left(  \sigma_{I,\alpha^{\prime},\beta^{\prime}}\left(  1\right)
,\sigma_{I,\alpha^{\prime},\beta^{\prime}}\left(  2\right)  ,\ldots
,\sigma_{I,\alpha^{\prime},\beta^{\prime}}\left(  n\right)  \right) \\
&  =\left(  a_{\alpha^{\prime}\left(  1\right)  },a_{\alpha^{\prime}\left(
2\right)  },\ldots,a_{\alpha^{\prime}\left(  k\right)  },b_{\beta^{\prime
}\left(  1\right)  },b_{\beta^{\prime}\left(  2\right)  },\ldots
,b_{\beta^{\prime}\left(  n-k\right)  }\right)  .
\end{align*}
Comparing this with%
\begin{align*}
&  \left(  \sigma_{I,\alpha^{\prime},\beta^{\prime}}\left(  1\right)
,\sigma_{I,\alpha^{\prime},\beta^{\prime}}\left(  2\right)  ,\ldots
,\sigma_{I,\alpha^{\prime},\beta^{\prime}}\left(  n\right)  \right) \\
&  =\left(  \sigma_{I,\alpha,\beta}\left(  1\right)  ,\sigma_{I,\alpha,\beta
}\left(  2\right)  ,\ldots,\sigma_{I,\alpha,\beta}\left(  n\right)  \right)
\ \ \ \ \ \ \ \ \ \ \left(  \text{since }\sigma_{I,\alpha^{\prime}%
,\beta^{\prime}}=\sigma_{I,\alpha,\beta}\right) \\
&  =\left(  a_{\alpha\left(  1\right)  },a_{\alpha\left(  2\right)  }%
,\ldots,a_{\alpha\left(  k\right)  },b_{\beta\left(  1\right)  }%
,b_{\beta\left(  2\right)  },\ldots,b_{\beta\left(  n-k\right)  }\right)  ,
\end{align*}
we obtain%
\begin{align*}
&  \left(  a_{\alpha\left(  1\right)  },a_{\alpha\left(  2\right)  }%
,\ldots,a_{\alpha\left(  k\right)  },b_{\beta\left(  1\right)  }%
,b_{\beta\left(  2\right)  },\ldots,b_{\beta\left(  n-k\right)  }\right) \\
&  =\left(  a_{\alpha^{\prime}\left(  1\right)  },a_{\alpha^{\prime}\left(
2\right)  },\ldots,a_{\alpha^{\prime}\left(  k\right)  },b_{\beta^{\prime
}\left(  1\right)  },b_{\beta^{\prime}\left(  2\right)  },\ldots
,b_{\beta^{\prime}\left(  n-k\right)  }\right)  .
\end{align*}
In other words,%
\begin{equation}
\left(  a_{\alpha\left(  i\right)  }=a_{\alpha^{\prime}\left(  i\right)
}\text{ for every }i\in\left\{  1,2,\ldots,k\right\}  \right)
\label{sol.Ialbe.c.short.inj.pf.1}%
\end{equation}
and%
\begin{equation}
\left(  b_{\beta\left(  j\right)  }=b_{\beta^{\prime}\left(  j\right)  }\text{
for every }j\in\left\{  1,2,\ldots,n-k\right\}  \right)  .
\label{sol.Ialbe.c.short.inj.pf.2}%
\end{equation}
\par
Now, we shall prove that $\alpha=\alpha^{\prime}$.
\par
Indeed, fix $i\in\left\{  1,2,\ldots,k\right\}  $. The list $\left(
a_{1},a_{2},\ldots,a_{k}\right)  $ has no repetitions. In other words, the
elements $a_{1},a_{2},\ldots,a_{k}$ are pairwise distinct. Thus, if $u$ and
$v$ are two elements of $\left\{  1,2,\ldots,k\right\}  $ such that
$a_{u}=a_{v}$, then $u=v$. Applying this to $u=\alpha\left(  i\right)  $ and
$v=\alpha^{\prime}\left(  i\right)  $, we obtain $\alpha\left(  i\right)
=\alpha^{\prime}\left(  i\right)  $ (since (\ref{sol.Ialbe.c.short.inj.pf.1})
yields $a_{\alpha\left(  i\right)  }=a_{\alpha^{\prime}\left(  i\right)  }$).
\par
Now, forget that we fixed $i$. We thus have proven that $\alpha\left(
i\right)  =\alpha^{\prime}\left(  i\right)  $ for every $i\in\left\{
1,2,\ldots,k\right\}  $. Thus, $\alpha=\alpha^{\prime}$ (since $\alpha$ and
$\alpha^{\prime}$ are permutations of $\left\{  1,2,\ldots,k\right\}  $).
\par
We have thus proven $\alpha=\alpha^{\prime}$ using the equalities
(\ref{sol.Ialbe.c.short.inj.pf.1}). Similarly, we can prove $\beta
=\beta^{\prime}$ using the equalities (\ref{sol.Ialbe.c.short.inj.pf.2}).
\par
Now, $\left(  \underbrace{\alpha}_{=\alpha^{\prime}},\underbrace{\beta
}_{=\beta^{\prime}}\right)  =\left(  \alpha^{\prime},\beta^{\prime}\right)  $.
\par
Now, forget that we fixed $\left(  \alpha,\beta\right)  $ and $\left(
\alpha^{\prime},\beta^{\prime}\right)  $. We thus have proven that if $\left(
\alpha,\beta\right)  $ and $\left(  \alpha^{\prime},\beta^{\prime}\right)  $
are two elements of $S_{k}\times S_{n-k}$ satisfying $\mu\left(  \alpha
,\beta\right)  =\mu\left(  \alpha^{\prime},\beta^{\prime}\right)  $, then
$\left(  \alpha,\beta\right)  =\left(  \alpha^{\prime},\beta^{\prime}\right)
$. In other words, the map $\mu$ is injective. Qed.}.
\end{vershort}

\begin{verlong}
For every $\left(  \alpha,\beta\right)  \in S_{k}\times S_{n-k}$ and every
$i\in\left\{  1,2,\ldots,k\right\}  $, we have%
\begin{equation}
a_{\alpha\left(  i\right)  }=\sigma_{I,\alpha,\beta}\left(  i\right)
\label{sol.Ialbe.c.aai=}%
\end{equation}
\footnote{\textit{Proof of (\ref{sol.Ialbe.c.aai=}):} Let $\left(
\alpha,\beta\right)  \in S_{k}\times S_{n-k}$.
\par
The permutation $\sigma_{I,\alpha,\beta}$ is the unique $\sigma\in S_{n}$
satisfying%
\begin{equation}
\left(  \sigma\left(  1\right)  ,\sigma\left(  2\right)  ,\ldots,\sigma\left(
n\right)  \right)  =\left(  a_{\alpha\left(  1\right)  },a_{\alpha\left(
2\right)  },\ldots,a_{\alpha\left(  k\right)  },b_{\beta\left(  1\right)
},b_{\beta\left(  2\right)  },\ldots,b_{\beta\left(  n-k\right)  }\right)  .
\label{sol.Ialbe.c.aai=.pf.1}%
\end{equation}
Thus, $\sigma_{I,\alpha,\beta}$ is a $\sigma\in S_{n}$ satisfying
(\ref{sol.Ialbe.c.aai=.pf.1}). In other words, we have%
\begin{align*}
&  \left(  \sigma_{I,\alpha,\beta}\left(  1\right)  ,\sigma_{I,\alpha,\beta
}\left(  2\right)  ,\ldots,\sigma_{I,\alpha,\beta}\left(  n\right)  \right) \\
&  =\left(  a_{\alpha\left(  1\right)  },a_{\alpha\left(  2\right)  }%
,\ldots,a_{\alpha\left(  k\right)  },b_{\beta\left(  1\right)  }%
,b_{\beta\left(  2\right)  },\ldots,b_{\beta\left(  n-k\right)  }\right)  .
\end{align*}
Now,%
\begin{align*}
&  \left(  \sigma_{I,\alpha,\beta}\left(  1\right)  ,\sigma_{I,\alpha,\beta
}\left(  2\right)  ,\ldots,\sigma_{I,\alpha,\beta}\left(  k\right)  \right) \\
&  =\left(  \text{the list of the first }k\text{ entries of the list
}\underbrace{\left(  \sigma_{I,\alpha,\beta}\left(  1\right)  ,\sigma
_{I,\alpha,\beta}\left(  2\right)  ,\ldots,\sigma_{I,\alpha,\beta}\left(
n\right)  \right)  }_{=\left(  a_{\alpha\left(  1\right)  },a_{\alpha\left(
2\right)  },\ldots,a_{\alpha\left(  k\right)  },b_{\beta\left(  1\right)
},b_{\beta\left(  2\right)  },\ldots,b_{\beta\left(  n-k\right)  }\right)
}\right) \\
&  =\left(  \text{the list of the first }k\text{ entries of the list }\left(
a_{\alpha\left(  1\right)  },a_{\alpha\left(  2\right)  },\ldots
,a_{\alpha\left(  k\right)  },b_{\beta\left(  1\right)  },b_{\beta\left(
2\right)  },\ldots,b_{\beta\left(  n-k\right)  }\right)  \right) \\
&  =\left(  a_{\alpha\left(  1\right)  },a_{\alpha\left(  2\right)  }%
,\ldots,a_{\alpha\left(  k\right)  }\right)  .
\end{align*}
In other words, $\left(  a_{\alpha\left(  1\right)  },a_{\alpha\left(
2\right)  },\ldots,a_{\alpha\left(  k\right)  }\right)  =\left(
\sigma_{I,\alpha,\beta}\left(  1\right)  ,\sigma_{I,\alpha,\beta}\left(
2\right)  ,\ldots,\sigma_{I,\alpha,\beta}\left(  k\right)  \right)  $. In
other words, $a_{\alpha\left(  i\right)  }=\sigma_{I,\alpha,\beta}\left(
i\right)  $ for every $i\in\left\{  1,2,\ldots,k\right\}  $. This proves
(\ref{sol.Ialbe.c.aai=}).}. For every $\left(  \alpha,\beta\right)  \in
S_{k}\times S_{n-k}$ and every $j\in\left\{  1,2,\ldots,n-k\right\}  $, we
have%
\begin{equation}
b_{\beta\left(  j\right)  }=\sigma_{I,\alpha,\beta}\left(  k+j\right)
\label{sol.Ialbe.c.bbj=}%
\end{equation}
\footnote{\textit{Proof of (\ref{sol.Ialbe.c.bbj=}):} Let $\left(
\alpha,\beta\right)  \in S_{k}\times S_{n-k}$.
\par
The permutation $\sigma_{I,\alpha,\beta}$ is the unique $\sigma\in S_{n}$
satisfying%
\begin{equation}
\left(  \sigma\left(  1\right)  ,\sigma\left(  2\right)  ,\ldots,\sigma\left(
n\right)  \right)  =\left(  a_{\alpha\left(  1\right)  },a_{\alpha\left(
2\right)  },\ldots,a_{\alpha\left(  k\right)  },b_{\beta\left(  1\right)
},b_{\beta\left(  2\right)  },\ldots,b_{\beta\left(  n-k\right)  }\right)  .
\label{sol.Ialbe.c.bbj=.pf.1}%
\end{equation}
Thus, $\sigma_{I,\alpha,\beta}$ is a $\sigma\in S_{n}$ satisfying
(\ref{sol.Ialbe.c.bbj=.pf.1}). In other words, we have%
\begin{align*}
&  \left(  \sigma_{I,\alpha,\beta}\left(  1\right)  ,\sigma_{I,\alpha,\beta
}\left(  2\right)  ,\ldots,\sigma_{I,\alpha,\beta}\left(  n\right)  \right) \\
&  =\left(  a_{\alpha\left(  1\right)  },a_{\alpha\left(  2\right)  }%
,\ldots,a_{\alpha\left(  k\right)  },b_{\beta\left(  1\right)  }%
,b_{\beta\left(  2\right)  },\ldots,b_{\beta\left(  n-k\right)  }\right)  .
\end{align*}
Now,%
\begin{align*}
&  \left(  \sigma_{I,\alpha,\beta}\left(  k+1\right)  ,\sigma_{I,\alpha,\beta
}\left(  k+2\right)  ,\ldots,\sigma_{I,\alpha,\beta}\left(  n\right)  \right)
\\
&  =\left(  \text{the list of the last }n-k\text{ entries of the list
}\underbrace{\left(  \sigma_{I,\alpha,\beta}\left(  1\right)  ,\sigma
_{I,\alpha,\beta}\left(  2\right)  ,\ldots,\sigma_{I,\alpha,\beta}\left(
n\right)  \right)  }_{=\left(  a_{\alpha\left(  1\right)  },a_{\alpha\left(
2\right)  },\ldots,a_{\alpha\left(  k\right)  },b_{\beta\left(  1\right)
},b_{\beta\left(  2\right)  },\ldots,b_{\beta\left(  n-k\right)  }\right)
}\right) \\
&  =\left(  \text{the list of the last }n-k\text{ entries of the list }\left(
a_{\alpha\left(  1\right)  },a_{\alpha\left(  2\right)  },\ldots
,a_{\alpha\left(  k\right)  },b_{\beta\left(  1\right)  },b_{\beta\left(
2\right)  },\ldots,b_{\beta\left(  n-k\right)  }\right)  \right) \\
&  =\left(  b_{\beta\left(  1\right)  },b_{\beta\left(  2\right)  }%
,\ldots,b_{\beta\left(  n-k\right)  }\right)  .
\end{align*}
In other words, $\left(  b_{\beta\left(  1\right)  },b_{\beta\left(  2\right)
},\ldots,b_{\beta\left(  n-k\right)  }\right)  =\left(  \sigma_{I,\alpha
,\beta}\left(  k+1\right)  ,\sigma_{I,\alpha,\beta}\left(  k+2\right)
,\ldots,\sigma_{I,\alpha,\beta}\left(  n\right)  \right)  $. In other words,
$b_{\beta\left(  j\right)  }=\sigma_{I,\alpha,\beta}\left(  k+j\right)  $ for
every $i\in\left\{  1,2,\ldots,k\right\}  $. This proves
(\ref{sol.Ialbe.c.bbj=}).}.

The map $\mu$ is injective\footnote{\textit{Proof.} Let $\eta$ and $\xi$ be
two elements of $S_{k}\times S_{n-k}$ such that $\mu\left(  \eta\right)
=\mu\left(  \xi\right)  $. We shall prove that $\eta=\xi$.
\par
We have $\eta\in S_{k}\times S_{n-k}$. Hence, $\eta=\left(  \alpha
,\beta\right)  $ for some $\alpha\in S_{k}$ and $\beta\in S_{n-k}$. Consider
these $\alpha$ and $\beta$.
\par
We have $\xi\in S_{k}\times S_{n-k}$. Hence, $\xi=\left(  \alpha^{\prime
},\beta^{\prime}\right)  $ for some $\alpha^{\prime}\in S_{k}$ and
$\beta^{\prime}\in S_{n-k}$. Consider these $\alpha^{\prime}$ and
$\beta^{\prime}$.
\par
We have $\eta=\left(  \alpha,\beta\right)  $ and thus $\mu\left(  \eta\right)
=\mu\left(  \alpha,\beta\right)  =\sigma_{I,\alpha,\beta}$ (by the definition
of the map $\mu$).
\par
We have $\xi=\left(  \alpha^{\prime},\beta^{\prime}\right)  $ and thus
$\mu\left(  \xi\right)  =\mu\left(  \alpha^{\prime},\beta^{\prime}\right)
=\sigma_{I,\alpha^{\prime},\beta^{\prime}}$ (by the definition of the map
$\xi$).
\par
From $\mu\left(  \eta\right)  =\sigma_{I,\alpha,\beta}$, we obtain
$\sigma_{I,\alpha,\beta}=\mu\left(  \eta\right)  =\mu\left(  \xi\right)
=\sigma_{I,\alpha^{\prime},\beta^{\prime}}$.
\par
We have $\alpha\in S_{k}$. In other words, $\alpha$ is a permutation of
$\left\{  1,2,\ldots,k\right\}  $ (since $S_{k}$ is the set of all
permutations of $\left\{  1,2,\ldots,k\right\}  $). In other words, $\alpha$
is a bijective map $\left\{  1,2,\ldots,k\right\}  \rightarrow\left\{
1,2,\ldots,k\right\}  $.
\par
Let $i\in\left\{  1,2,\ldots,k\right\}  $. Then, $\alpha\left(  i\right)  $
and $\alpha^{\prime}\left(  i\right)  $ are elements of $\left\{
1,2,\ldots,k\right\}  $.
\par
Furthermore, (\ref{sol.Ialbe.c.aai=}) shows that $a_{\alpha\left(  i\right)
}=\underbrace{\sigma_{I,\alpha,\beta}}_{=\sigma_{I,\alpha^{\prime}%
,\beta^{\prime}}}\left(  i\right)  =\sigma_{I,\alpha^{\prime},\beta^{\prime}%
}\left(  i\right)  $. But (\ref{sol.Ialbe.c.aai=}) (applied to $\left(
\alpha^{\prime},\beta^{\prime}\right)  $ instead of $\left(  \alpha
,\beta\right)  $) yields $a_{\alpha^{\prime}\left(  i\right)  }=\sigma
_{I,\alpha^{\prime},\beta^{\prime}}\left(  i\right)  $. Comparing this with
$a_{\alpha\left(  i\right)  }=\sigma_{I,\alpha^{\prime},\beta^{\prime}}\left(
i\right)  $, we obtain $a_{\alpha\left(  i\right)  }=a_{\alpha^{\prime}\left(
i\right)  }$.
\par
We know that $\left(  a_{1},a_{2},\ldots,a_{k}\right)  $ is a list with no
repetitions. In other words, the elements $a_{1},a_{2},\ldots,a_{k}$ are
pairwise distinct. In other words, if $u$ and $v$ are two elements of
$\left\{  1,2,\ldots,k\right\}  $ such that $a_{u}=a_{v}$, then $u=v$.
Applying this to $u=\alpha\left(  i\right)  $ and $v=\alpha^{\prime}\left(
i\right)  $, we obtain $\alpha\left(  i\right)  =\alpha^{\prime}\left(
i\right)  $ (since $a_{\alpha\left(  i\right)  }=a_{\alpha^{\prime}\left(
i\right)  }$).
\par
Now, forget that we fixed $i$. We thus have shown that $\alpha\left(
i\right)  =\alpha^{\prime}\left(  i\right)  $ for every $i\in\left\{
1,2,\ldots,k\right\}  $. In other words, $\alpha=\alpha^{\prime}$ (since
$\alpha$ and $\alpha^{\prime}$ are maps $\left\{  1,2,\ldots,k\right\}
\rightarrow\left\{  1,2,\ldots,k\right\}  $).
\par
We have $\beta\in S_{n-k}$. In other words, $\beta$ is a permutation of
$\left\{  1,2,\ldots,n-k\right\}  $ (since $S_{n-k}$ is the set of all
permutations of $\left\{  1,2,\ldots,n-k\right\}  $). In other words, $\beta$
is a bijective map $\left\{  1,2,\ldots,n-k\right\}  \rightarrow\left\{
1,2,\ldots,n-k\right\}  $.
\par
Let $j\in\left\{  1,2,\ldots,n-k\right\}  $. Then, $\beta\left(  j\right)  $
and $\beta^{\prime}\left(  j\right)  $ are elements of $\left\{
1,2,\ldots,n-k\right\}  $.
\par
Furthermore, (\ref{sol.Ialbe.c.bbj=}) shows that $b_{\beta\left(  j\right)
}=\underbrace{\sigma_{I,\alpha,\beta}}_{=\sigma_{I,\alpha^{\prime}%
,\beta^{\prime}}}\left(  k+j\right)  =\sigma_{I,\alpha^{\prime},\beta^{\prime
}}\left(  k+j\right)  $. But (\ref{sol.Ialbe.c.bbj=}) (applied to $\left(
\alpha^{\prime},\beta^{\prime}\right)  $ instead of $\left(  \alpha
,\beta\right)  $) yields $b_{\beta^{\prime}\left(  j\right)  }=\sigma
_{I,\alpha^{\prime},\beta^{\prime}}\left(  k+j\right)  $. Comparing this with
$b_{\beta\left(  j\right)  }=\sigma_{I,\alpha^{\prime},\beta^{\prime}}\left(
k+j\right)  $, we obtain $b_{\beta\left(  j\right)  }=b_{\beta^{\prime}\left(
j\right)  }$.
\par
We know that $\left(  b_{1},b_{2},\ldots,b_{n-k}\right)  $ is a list with no
repetitions. In other words, the elements $b_{1},b_{2},\ldots,b_{n-k}$ are
pairwise distinct. In other words, if $u$ and $v$ are two elements of
$\left\{  1,2,\ldots,n-k\right\}  $ such that $b_{u}=b_{v}$, then $u=v$.
Applying this to $u=\beta\left(  j\right)  $ and $v=\beta^{\prime}\left(
j\right)  $, we obtain $\beta\left(  j\right)  =\beta^{\prime}\left(
j\right)  $ (since $b_{\beta\left(  j\right)  }=b_{\beta^{\prime}\left(
j\right)  }$).
\par
Now, forget that we fixed $j$. We thus have shown that $\beta\left(  j\right)
=\beta^{\prime}\left(  j\right)  $ for every $j\in\left\{  1,2,\ldots
,n-k\right\}  $. In other words, $\beta=\beta^{\prime}$ (since $\beta$ and
$\beta^{\prime}$ are maps $\left\{  1,2,\ldots,n-k\right\}  \rightarrow
\left\{  1,2,\ldots,n-k\right\}  $).
\par
Now, we have shown that $\alpha=\alpha^{\prime}$ and $\beta=\beta^{\prime}$.
Now, $\eta=\left(  \underbrace{\alpha}_{=\alpha^{\prime}},\underbrace{\beta
}_{=\beta^{\prime}}\right)  =\left(  \alpha^{\prime},\beta^{\prime}\right)
=\xi$ (since $\xi=\left(  \alpha^{\prime},\beta^{\prime}\right)  $).
\par
Now, forget that we fixed $\eta$ and $\xi$. We thus have proven that if $\eta$
and $\xi$ are two elements of $S_{k}\times S_{n-k}$ such that $\mu\left(
\eta\right)  =\mu\left(  \xi\right)  $, then $\eta=\xi$. In other words, the
map $\mu$ is injective. Qed.}.
\end{verlong}

Our next goal is to show that the map $\mu$ is surjective.

In fact, let $\gamma\in\left\{  \tau\in S_{n}\ \mid\ \tau\left(  \left\{
1,2,\ldots,k\right\}  \right)  =I\right\}  $. We shall construct an $\left(
\alpha,\beta\right)  \in S_{k}\times S_{n-k}$ such that $\mu\left(
\alpha,\beta\right)  =\gamma$.

\begin{vershort}
We have $\gamma\in\left\{  \tau\in S_{n}\ \mid\ \tau\left(  \left\{
1,2,\ldots,k\right\}  \right)  =I\right\}  $. In other words, $\gamma$ is an
element of $S_{n}$ and satisfies $\gamma\left(  \left\{  1,2,\ldots,k\right\}
\right)  =I$.
\end{vershort}

\begin{verlong}
We have $\gamma\in\left\{  \tau\in S_{n}\ \mid\ \tau\left(  \left\{
1,2,\ldots,k\right\}  \right)  =I\right\}  $. In other words, $\gamma$ is an
element $\tau$ of $S_{n}$ satisfying $\tau\left(  \left\{  1,2,\ldots
,k\right\}  \right)  =I$. In other words, $\gamma$ is an element of $S_{n}$
and satisfies $\gamma\left(  \left\{  1,2,\ldots,k\right\}  \right)  =I$.
\end{verlong}

\begin{vershort}
We have $\gamma\in S_{n}$. In other words, $\gamma$ is a permutation of
$\left\{  1,2,\ldots,n\right\}  $. Thus, $\gamma$ is a bijective map $\left\{
1,2,\ldots,n\right\}  \rightarrow\left\{  1,2,\ldots,n\right\}  $, and
therefore also an injective map.
\end{vershort}

\begin{verlong}
We have $\gamma\in S_{n}$. In other words, $\gamma$ is a permutation of
$\left\{  1,2,\ldots,n\right\}  $ (since $S_{n}$ is the set of all
permutations of $\left\{  1,2,\ldots,n\right\}  $). In other words, $\gamma$
is a bijective map $\left\{  1,2,\ldots,n\right\}  \rightarrow\left\{
1,2,\ldots,n\right\}  $. The map $\gamma$ is bijective and thus injective.
\end{verlong}

Every $i\in\left\{  1,2,\ldots,k\right\}  $ satisfies $\gamma\left(  i\right)
\in I$\ \ \ \ \footnote{\textit{Proof.} Let $i\in\left\{  1,2,\ldots
,k\right\}  $. Then, $\gamma\left(  \underbrace{i}_{\in\left\{  1,2,\ldots
,k\right\}  }\right)  \in\gamma\left(  \left\{  1,2,\ldots,k\right\}  \right)
=I$, qed.}. Thus, $\gamma\left(  1\right)  ,\gamma\left(  2\right)
,\ldots,\gamma\left(  k\right)  $ are $k$ elements of $I$.

\begin{vershort}
If $u$ and $v$ are two distinct elements of $\left\{  1,2,\ldots,k\right\}  $,
then $\gamma\left(  u\right)  \neq\gamma\left(  v\right)  $ (since the map
$\gamma$ is injective). In other words, the $k$ elements $\gamma\left(
1\right)  ,\gamma\left(  2\right)  ,\ldots,\gamma\left(  k\right)  $ are
pairwise distinct.
\end{vershort}

\begin{verlong}
Now, if $u$ and $v$ are two distinct elements of $\left\{  1,2,\ldots
,k\right\}  $, then $\gamma\left(  u\right)  \neq\gamma\left(  v\right)
$\ \ \ \ \footnote{\textit{Proof.} Let $u$ and $v$ be two distinct elements of
$\left\{  1,2,\ldots,k\right\}  $. Then, $u\in\left\{  1,2,\ldots,k\right\}
\subseteq\left\{  1,2,\ldots,n\right\}  $ (since $k\leq n$) and $v\in\left\{
1,2,\ldots,k\right\}  \subseteq\left\{  1,2,\ldots,n\right\}  $. Now, $u\neq
v$ (since $u$ and $v$ are distinct), so that $\gamma\left(  u\right)
\neq\gamma\left(  v\right)  $ (since the map $\gamma$ is injective). Qed.}. In
other words, the elements $\gamma\left(  1\right)  ,\gamma\left(  2\right)
,\ldots,\gamma\left(  k\right)  $ are pairwise distinct. Hence, $\gamma\left(
1\right)  ,\gamma\left(  2\right)  ,\ldots,\gamma\left(  k\right)  $ are $k$
pairwise distinct elements of $I$ (since $\gamma\left(  1\right)
,\gamma\left(  2\right)  ,\ldots,\gamma\left(  k\right)  $ are $k$ elements of
$I$).
\end{verlong}

On the other hand, recall that $\left(  a_{1},a_{2},\ldots,a_{k}\right)  $ is
a list of all elements of $I$ (with no repetitions). Thus, Lemma
\ref{lem.sol.Ialbe.inclist2} (applied to $I$, $k$, $\left(  a_{1},a_{2}%
,\ldots,a_{k}\right)  $ and \newline$\left(  \gamma\left(  1\right)
,\gamma\left(  2\right)  ,\ldots,\gamma\left(  k\right)  \right)  $ instead of
$S$, $s$, $\left(  c_{1},c_{2},\ldots,c_{s}\right)  $ and $\left(  p_{1}%
,p_{2},\ldots,p_{s}\right)  $) yields that there exists a $\pi\in S_{k}$ such
that $\left(  \gamma\left(  1\right)  ,\gamma\left(  2\right)  ,\ldots
,\gamma\left(  k\right)  \right)  =\left(  a_{\pi\left(  1\right)  }%
,a_{\pi\left(  2\right)  },\ldots,a_{\pi\left(  k\right)  }\right)  $. Denote
this $\pi$ by $\alpha$. Thus, $\alpha$ is a $\pi\in S_{k}$ such that $\left(
\gamma\left(  1\right)  ,\gamma\left(  2\right)  ,\ldots,\gamma\left(
k\right)  \right)  =\left(  a_{\pi\left(  1\right)  },a_{\pi\left(  2\right)
},\ldots,a_{\pi\left(  k\right)  }\right)  $. In other words, $\alpha$ is an
element of $S_{k}$ and satisfies
\begin{equation}
\left(  \gamma\left(  1\right)  ,\gamma\left(  2\right)  ,\ldots,\gamma\left(
k\right)  \right)  =\left(  a_{\alpha\left(  1\right)  },a_{\alpha\left(
2\right)  },\ldots,a_{\alpha\left(  k\right)  }\right)  .
\label{sol.Ialbe.c.gamma1}%
\end{equation}

\begin{vershort}
Every $j\in\left\{  1,2,\ldots,n-k\right\}  $ satisfies $\gamma\left(
k+j\right)  \in\left\{  1,2,\ldots,n\right\}  \setminus I$%
\ \ \ \ \footnote{\textit{Proof.} Let $j\in\left\{  1,2,\ldots,n-k\right\}  $.
Then, $k+j\in\left\{  k+1,k+2,\ldots,n\right\}  \subseteq\left\{
1,2,\ldots,n\right\}  $. Hence, $\gamma\left(  k+j\right)  $ is well-defined.
\par
We must show that $\gamma\left(  k+j\right)  \in\left\{  1,2,\ldots,n\right\}
\setminus I$.
\par
Indeed, assume the contrary. Thus, we don't have $\gamma\left(  k+j\right)
\in\left\{  1,2,\ldots,n\right\}  \setminus I$. In other words, we have
$\gamma\left(  k+j\right)  \notin\left\{  1,2,\ldots,n\right\}  \setminus I$.
Combining $\gamma\left(  k+j\right)  \in\left\{  1,2,\ldots,n\right\}  $ with
$\gamma\left(  k+j\right)  \notin\left\{  1,2,\ldots,n\right\}  \setminus I$,
we obtain
\[
\gamma\left(  k+j\right)  \in\left\{  1,2,\ldots,n\right\}  \setminus\left(
\left\{  1,2,\ldots,n\right\}  \setminus I\right)  \subseteq I=\gamma\left(
\left\{  1,2,\ldots,k\right\}  \right)
\]
(since $\gamma\left(  \left\{  1,2,\ldots,k\right\}  \right)  =I$). In other
words, there exists some $v\in\left\{  1,2,\ldots,k\right\}  $ such that
$\gamma\left(  k+j\right)  =\gamma\left(  v\right)  $. Consider this $v$.
Since $\gamma\left(  k+j\right)  =\gamma\left(  v\right)  $, we have $k+j=v$
(because $\gamma$ is injective). Thus, $k+j=v\in\left\{  1,2,\ldots,k\right\}
$, so that $k+j\leq k$. Therefore, $j\leq k-k=0$. This contradicts
$j\in\left\{  1,2,\ldots,n-k\right\}  $. This contradiction shows that our
assumption was wrong. Hence, $\gamma\left(  k+j\right)  \in\left\{
1,2,\ldots,n\right\}  \setminus I$ is proven. Qed.}. Thus, $\gamma\left(
k+1\right)  ,\gamma\left(  k+2\right)  ,\ldots,\gamma\left(  n\right)  $ are
$n-k$ elements of $\left\{  1,2,\ldots,n\right\}  \setminus I$.
\end{vershort}

\begin{verlong}
Every $j\in\left\{  1,2,\ldots,n-k\right\}  $ satisfies $\gamma\left(
k+j\right)  \in\left\{  1,2,\ldots,n\right\}  \setminus I$%
\ \ \ \ \footnote{\textit{Proof.} Let $j\in\left\{  1,2,\ldots,n-k\right\}  $.
Then, $k+j\in\left\{  k+1,k+2,\ldots,n\right\}  \subseteq\left\{
1,2,\ldots,n\right\}  $. Hence, $\gamma\left(  k+j\right)  $ is well-defined.
\par
We must show that $\gamma\left(  k+j\right)  \in\left\{  1,2,\ldots,n\right\}
\setminus I$.
\par
Indeed, assume the contrary. Thus, we don't have $\gamma\left(  k+j\right)
\in\left\{  1,2,\ldots,n\right\}  \setminus I$. In other words, we have
$\gamma\left(  k+j\right)  \notin\left\{  1,2,\ldots,n\right\}  \setminus I$.
Combining $\gamma\left(  k+j\right)  \in\left\{  1,2,\ldots,n\right\}  $
(since $\gamma$ is a map $\left\{  1,2,\ldots,n\right\}  \rightarrow\left\{
1,2,\ldots,n\right\}  $) with $\gamma\left(  k+j\right)  \notin\left\{
1,2,\ldots,n\right\}  \setminus I$, we obtain
\begin{align*}
\gamma\left(  k+j\right)   &  \in\left\{  1,2,\ldots,n\right\}  \setminus
\left(  \left\{  1,2,\ldots,n\right\}  \setminus I\right)
=I\ \ \ \ \ \ \ \ \ \ \left(  \text{since }I\subseteq\left\{  1,2,\ldots
,n\right\}  \right) \\
&  =\left\{  a_{1},a_{2},\ldots,a_{k}\right\}  \ \ \ \ \ \ \ \ \ \ \left(
\text{since }\left(  a_{1},a_{2},\ldots,a_{k}\right)  \text{ is a list of all
elements of }I\right)  .
\end{align*}
In other words, $\gamma\left(  k+j\right)  =a_{p}$ for some $p\in\left\{
1,2,\ldots,k\right\}  $. Consider this $p$. Clearly, $\alpha^{-1}\left(
p\right)  $ is a well-defined element of $\left\{  1,2,\ldots,k\right\}  $
(since $\alpha$ is a permutation of $\left\{  1,2,\ldots,k\right\}  $, while
$p$ is an element of $\left\{  1,2,\ldots,k\right\}  $).
\par
But recall that $\left(  \gamma\left(  1\right)  ,\gamma\left(  2\right)
,\ldots,\gamma\left(  k\right)  \right)  =\left(  a_{\alpha\left(  1\right)
},a_{\alpha\left(  2\right)  },\ldots,a_{\alpha\left(  k\right)  }\right)  $.
In other words, $\gamma\left(  u\right)  =a_{\alpha\left(  u\right)  }$ for
every $u\in\left\{  1,2,\ldots,k\right\}  $. Applying this to $u=\alpha
^{-1}\left(  p\right)  $, we obtain $\gamma\left(  \alpha^{-1}\left(
p\right)  \right)  =a_{\alpha\left(  \alpha^{-1}\left(  p\right)  \right)
}=a_{p}$ (since $\alpha\left(  \alpha^{-1}\left(  p\right)  \right)  =p$).
Comparing this with $\gamma\left(  k+j\right)  =a_{p}$, we obtain
$\gamma\left(  k+j\right)  =\gamma\left(  \alpha^{-1}\left(  p\right)
\right)  $.
\par
But the map $\gamma$ is injective. In other words, if $u$ and $v$ are two
elements of $\left\{  1,2,\ldots,n\right\}  $ such that $\gamma\left(
u\right)  =\gamma\left(  v\right)  $, then $u=v$. Applying this to $u=k+j$ and
$v=\alpha^{-1}\left(  p\right)  $, we obtain $k+j=\alpha^{-1}\left(  p\right)
$ (since $\gamma\left(  k+j\right)  =\gamma\left(  \alpha^{-1}\left(
p\right)  \right)  $). But $\alpha^{-1}\left(  p\right)  \in\left\{
1,2,\ldots,k\right\}  $, so that $\alpha^{-1}\left(  p\right)  \leq k$. Hence,
$k+j=\alpha^{-1}\left(  p\right)  \leq k$, so that $j\leq k-k=0<1$. This
contradicts $j\geq1$ (since $j\in\left\{  1,2,\ldots,n-k\right\}  $). This
contradiction shows that our assumption was wrong. Hence, $\gamma\left(
k+j\right)  \in\left\{  1,2,\ldots,n\right\}  \setminus I$ is proven. Qed.}.
Thus, $\gamma\left(  k+1\right)  ,\gamma\left(  k+2\right)  ,\ldots
,\gamma\left(  k+\left(  n-k\right)  \right)  $ are $n-k$ elements of
$\left\{  1,2,\ldots,n\right\}  \setminus I$.
\end{verlong}

\begin{vershort}
If $u$ and $v$ are two distinct elements of $\left\{  k+1,k+2,\ldots
,n\right\}  $, then $\gamma\left(  u\right)  \neq\gamma\left(  v\right)  $
(since the map $\gamma$ is injective). In other words, the $n-k$ elements
\newline$\gamma\left(  k+1\right)  ,\gamma\left(  k+2\right)  ,\ldots
,\gamma\left(  n\right)  $ are pairwise distinct.
\end{vershort}

\begin{verlong}
Now, if $u$ and $v$ are two distinct elements of $\left\{  1,2,\ldots
,n-k\right\}  $, then $\gamma\left(  k+u\right)  \neq\gamma\left(  k+v\right)
$\ \ \ \ \footnote{\textit{Proof.} Let $u$ and $v$ be two distinct elements of
$\left\{  1,2,\ldots,n-k\right\}  $. Then, $u\in\left\{  1,2,\ldots
,n-k\right\}  $, so that $k+u\in\left\{  k+1,k+2,\ldots,n\right\}
\subseteq\left\{  1,2,\ldots,n\right\}  $.The same argument (applied to $v$
instead of $u$) shows that $k+v\in\left\{  1,2,\ldots,n\right\}  $. Now,
$u\neq v$ (since $u$ and $v$ are distinct), so that $k+u\neq k+v$. Hence,
$\gamma\left(  k+u\right)  \neq\gamma\left(  k+v\right)  $ (since the map
$\gamma$ is injective). Qed.}. In other words, the elements $\gamma\left(
k+1\right)  ,\gamma\left(  k+2\right)  ,\ldots,\gamma\left(  k+\left(
n-k\right)  \right)  $ are pairwise distinct. Hence, $\gamma\left(
k+1\right)  ,\gamma\left(  k+2\right)  ,\ldots,\gamma\left(  k+\left(
n-k\right)  \right)  $ are $n-k$ pairwise distinct elements of $\left\{
1,2,\ldots,n\right\}  \setminus I$ (since $\gamma\left(  k+1\right)
,\gamma\left(  k+2\right)  ,\ldots,\gamma\left(  k+\left(  n-k\right)
\right)  $ are $n-k$ elements of $\left\{  1,2,\ldots,n\right\}  \setminus I$).
\end{verlong}

\begin{vershort}
On the other hand, recall that $\left(  b_{1},b_{2},\ldots,b_{n-k}\right)  $
is a list of all elements of $\left\{  1,2,\ldots,n\right\}  \setminus I$
(with no repetitions). Thus, Lemma \ref{lem.sol.Ialbe.inclist2} (applied to
$\left\{  1,2,\ldots,n\right\}  \setminus I$, $n-k$, \newline$\left(
b_{1},b_{2},\ldots,b_{n-k}\right)  $ and $\left(  \gamma\left(  k+1\right)
,\gamma\left(  k+2\right)  ,\ldots,\gamma\left(  n\right)  \right)  $ instead
of $S$, $s$, $\left(  c_{1},c_{2},\ldots,c_{s}\right)  $ and $\left(
p_{1},p_{2},\ldots,p_{s}\right)  $) yields that there exists a $\pi\in
S_{n-k}$ such that \newline$\left(  \gamma\left(  k+1\right)  ,\gamma\left(
k+2\right)  ,\ldots,\gamma\left(  n\right)  \right)  =\left(  b_{\pi\left(
1\right)  },b_{\pi\left(  2\right)  },\ldots,b_{\pi\left(  n-k\right)
}\right)  $. Denote this $\pi$ by $\beta$. Thus, $\beta$ is an element of
$S_{n-k}$ and satisfies
\begin{equation}
\left(  \gamma\left(  k+1\right)  ,\gamma\left(  k+2\right)  ,\ldots
,\gamma\left(  n\right)  \right)  =\left(  b_{\beta\left(  1\right)
},b_{\beta\left(  2\right)  },\ldots,b_{\beta\left(  n-k\right)  }\right)  .
\label{sol.Ialbe.c.short.gamma2}%
\end{equation}

\end{vershort}

\begin{verlong}
On the other hand, recall that $\left(  b_{1},b_{2},\ldots,b_{n-k}\right)  $
is a list of all elements of $\left\{  1,2,\ldots,n\right\}  \setminus I$
(with no repetitions). Thus, Lemma \ref{lem.sol.Ialbe.inclist2} (applied to
$\left\{  1,2,\ldots,n\right\}  \setminus I$, $n-k$, $\left(  b_{1}%
,b_{2},\ldots,b_{n-k}\right)  $ and $\left(  \gamma\left(  k+1\right)
,\gamma\left(  k+2\right)  ,\ldots,\gamma\left(  k+\left(  n-k\right)
\right)  \right)  $ instead of $S$, $s$, \newline$\left(  c_{1},c_{2}%
,\ldots,c_{s}\right)  $ and $\left(  p_{1},p_{2},\ldots,p_{s}\right)  $)
yields that there exists a $\pi\in S_{n-k}$ such that $\left(  \gamma\left(
k+1\right)  ,\gamma\left(  k+2\right)  ,\ldots,\gamma\left(  k+\left(
n-k\right)  \right)  \right)  =\left(  b_{\pi\left(  1\right)  },b_{\pi\left(
2\right)  },\ldots,b_{\pi\left(  n-k\right)  }\right)  $. Denote this $\pi$ by
$\beta$. Thus, $\beta$ is a $\pi\in S_{n-k}$ such that $\left(  \gamma\left(
k+1\right)  ,\gamma\left(  k+2\right)  ,\ldots,\gamma\left(  k+\left(
n-k\right)  \right)  \right)  =\left(  b_{\pi\left(  1\right)  },b_{\pi\left(
2\right)  },\ldots,b_{\pi\left(  n-k\right)  }\right)  $. In other words,
$\beta$ is an element of $S_{n-k}$ and satisfies
\begin{equation}
\left(  \gamma\left(  k+1\right)  ,\gamma\left(  k+2\right)  ,\ldots
,\gamma\left(  k+\left(  n-k\right)  \right)  \right)  =\left(  b_{\beta
\left(  1\right)  },b_{\beta\left(  2\right)  },\ldots,b_{\beta\left(
n-k\right)  }\right)  . \label{sol.Ialbe.c.gamma2}%
\end{equation}

\end{verlong}

\begin{verlong}
The permutation $\sigma_{I,\alpha,\beta}$ is the unique $\sigma\in S_{n}$
satisfying%
\begin{equation}
\left(  \sigma\left(  1\right)  ,\sigma\left(  2\right)  ,\ldots,\sigma\left(
n\right)  \right)  =\left(  a_{\alpha\left(  1\right)  },a_{\alpha\left(
2\right)  },\ldots,a_{\alpha\left(  k\right)  },b_{\beta\left(  1\right)
},b_{\beta\left(  2\right)  },\ldots,b_{\beta\left(  n-k\right)  }\right)  .
\label{sol.Ialbe.c.sigma1}%
\end{equation}
Thus, $\sigma_{I,\alpha,\beta}$ is a $\sigma\in S_{n}$ satisfying
(\ref{sol.Ialbe.c.sigma1}). In other words, we have%
\begin{align}
&  \left(  \sigma_{I,\alpha,\beta}\left(  1\right)  ,\sigma_{I,\alpha,\beta
}\left(  2\right)  ,\ldots,\sigma_{I,\alpha,\beta}\left(  n\right)  \right)
\nonumber\\
&  =\left(  a_{\alpha\left(  1\right)  },a_{\alpha\left(  2\right)  }%
,\ldots,a_{\alpha\left(  k\right)  },b_{\beta\left(  1\right)  }%
,b_{\beta\left(  2\right)  },\ldots,b_{\beta\left(  n-k\right)  }\right)  .
\label{sol.Ialbe.c.sigma2}%
\end{align}

\end{verlong}

\begin{vershort}
Now, let us introduce a notation: If $\left(  x_{1},x_{2},\ldots,x_{u}\right)
$ and $\left(  y_{1},y_{2},\ldots,y_{v}\right)  $ are two lists, then we shall
let $\left(  x_{1},x_{2},\ldots,x_{u}\right)  \ast\left(  y_{1},y_{2}%
,\ldots,y_{v}\right)  $ denote the list \newline$\left(  x_{1},x_{2}%
,\ldots,x_{u},y_{1},y_{2},\ldots,y_{v}\right)  $. Then,%
\begin{align*}
&  \left(  \gamma\left(  1\right)  ,\gamma\left(  2\right)  ,\ldots
,\gamma\left(  n\right)  \right) \\
&  =\underbrace{\left(  \gamma\left(  1\right)  ,\gamma\left(  2\right)
,\ldots,\gamma\left(  k\right)  \right)  }_{\substack{=\left(  a_{\alpha
\left(  1\right)  },a_{\alpha\left(  2\right)  },\ldots,a_{\alpha\left(
k\right)  }\right)  \\\text{(by (\ref{sol.Ialbe.c.gamma1}))}}}\ast
\underbrace{\left(  \gamma\left(  k+1\right)  ,\gamma\left(  k+2\right)
,\ldots,\gamma\left(  n\right)  \right)  }_{\substack{=\left(  b_{\beta\left(
1\right)  },b_{\beta\left(  2\right)  },\ldots,b_{\beta\left(  n-k\right)
}\right)  \\\text{(by (\ref{sol.Ialbe.c.short.gamma2}))}}}\\
&  =\left(  a_{\alpha\left(  1\right)  },a_{\alpha\left(  2\right)  }%
,\ldots,a_{\alpha\left(  k\right)  }\right)  \ast\left(  b_{\beta\left(
1\right)  },b_{\beta\left(  2\right)  },\ldots,b_{\beta\left(  n-k\right)
}\right) \\
&  =\left(  a_{\alpha\left(  1\right)  },a_{\alpha\left(  2\right)  }%
,\ldots,a_{\alpha\left(  k\right)  },b_{\beta\left(  1\right)  }%
,b_{\beta\left(  2\right)  },\ldots,b_{\beta\left(  n-k\right)  }\right) \\
&  =\left(  \sigma_{I,\alpha,\beta}\left(  1\right)  ,\sigma_{I,\alpha,\beta
}\left(  2\right)  ,\ldots,\sigma_{I,\alpha,\beta}\left(  n\right)  \right)
\ \ \ \ \ \ \ \ \ \ \left(  \text{by (\ref{sol.Ialbe.c.short.siIAB})}\right)
.
\end{align*}
In other words, $\gamma\left(  i\right)  =\sigma_{I,\alpha,\beta}\left(
i\right)  $ for every $i\in\left\{  1,2,\ldots,n\right\}  $. In other words,
$\gamma=\sigma_{I,\alpha,\beta}$ (since $\gamma$ and $\sigma_{I,\alpha,\beta}$
are two maps $\left\{  1,2,\ldots,n\right\}  \rightarrow\left\{
1,2,\ldots,n\right\}  $). But the definition of $\mu$ yields $\mu\left(
\alpha,\beta\right)  =\sigma_{I,\alpha,\beta}$. Comparing this with
$\gamma=\sigma_{I,\alpha,\beta}$, we obtain
\[
\gamma=\mu\underbrace{\left(  \alpha,\beta\right)  }_{\in S_{k}\times S_{n-k}%
}\in\mu\left(  S_{k}\times S_{n-k}\right)  .
\]

\end{vershort}

\begin{verlong}
Now, let us introduce a notation: If $\left(  x_{1},x_{2},\ldots,x_{u}\right)
$ and $\left(  y_{1},y_{2},\ldots,y_{v}\right)  $ are two lists, then we shall
let $\left(  x_{1},x_{2},\ldots,x_{u}\right)  \ast\left(  y_{1},y_{2}%
,\ldots,y_{v}\right)  $ denote the list \newline$\left(  x_{1},x_{2}%
,\ldots,x_{u},y_{1},y_{2},\ldots,y_{v}\right)  $. Then,%
\begin{align*}
&  \left(  \gamma\left(  1\right)  ,\gamma\left(  2\right)  ,\ldots
,\gamma\left(  n\right)  \right) \\
&  =\left(  \gamma\left(  1\right)  ,\gamma\left(  2\right)  ,\ldots
,\gamma\left(  k\right)  \right)  \ast\left(  \gamma\left(  k+1\right)
,\gamma\left(  k+2\right)  ,\ldots,\gamma\left(  \underbrace{n}_{=k+\left(
n-k\right)  }\right)  \right) \\
&  =\underbrace{\left(  \gamma\left(  1\right)  ,\gamma\left(  2\right)
,\ldots,\gamma\left(  k\right)  \right)  }_{\substack{=\left(  a_{\alpha
\left(  1\right)  },a_{\alpha\left(  2\right)  },\ldots,a_{\alpha\left(
k\right)  }\right)  \\\text{(by (\ref{sol.Ialbe.c.gamma1}))}}}\ast
\underbrace{\left(  \gamma\left(  k+1\right)  ,\gamma\left(  k+2\right)
,\ldots,\gamma\left(  k+\left(  n-k\right)  \right)  \right)  }%
_{\substack{=\left(  b_{\beta\left(  1\right)  },b_{\beta\left(  2\right)
},\ldots,b_{\beta\left(  n-k\right)  }\right)  \\\text{(by
(\ref{sol.Ialbe.c.gamma2}))}}}\\
&  =\left(  a_{\alpha\left(  1\right)  },a_{\alpha\left(  2\right)  }%
,\ldots,a_{\alpha\left(  k\right)  }\right)  \ast\left(  b_{\beta\left(
1\right)  },b_{\beta\left(  2\right)  },\ldots,b_{\beta\left(  n-k\right)
}\right) \\
&  =\left(  a_{\alpha\left(  1\right)  },a_{\alpha\left(  2\right)  }%
,\ldots,a_{\alpha\left(  k\right)  },b_{\beta\left(  1\right)  }%
,b_{\beta\left(  2\right)  },\ldots,b_{\beta\left(  n-k\right)  }\right) \\
&  =\left(  \sigma_{I,\alpha,\beta}\left(  1\right)  ,\sigma_{I,\alpha,\beta
}\left(  2\right)  ,\ldots,\sigma_{I,\alpha,\beta}\left(  n\right)  \right)
\ \ \ \ \ \ \ \ \ \ \left(  \text{by (\ref{sol.Ialbe.c.sigma2})}\right)  .
\end{align*}
In other words, $\gamma\left(  i\right)  =\sigma_{I,\alpha,\beta}\left(
i\right)  $ for every $i\in\left\{  1,2,\ldots,n\right\}  $. In other words,
$\gamma=\sigma_{I,\alpha,\beta}$ (since $\gamma$ and $\sigma_{I,\alpha,\beta}$
are two maps $\left\{  1,2,\ldots,n\right\}  \rightarrow\left\{
1,2,\ldots,n\right\}  $). But the definition of $\mu$ yields $\mu\left(
\alpha,\beta\right)  =\sigma_{I,\alpha,\beta}$. Comparing this with
$\gamma=\sigma_{I,\alpha,\beta}$, we obtain
\[
\gamma=\mu\underbrace{\left(  \alpha,\beta\right)  }_{\in S_{k}\times S_{n-k}%
}\in\mu\left(  S_{k}\times S_{n-k}\right)  .
\]

\end{verlong}

Now, forget that we fixed $\gamma$. We thus have shown that every
\newline$\gamma\in\left\{  \tau\in S_{n}\ \mid\ \tau\left(  \left\{
1,2,\ldots,k\right\}  \right)  =I\right\}  $ satisfies $\gamma\in\mu\left(
S_{k}\times S_{n-k}\right)  $. In other words,%
\[
\left\{  \tau\in S_{n}\ \mid\ \tau\left(  \left\{  1,2,\ldots,k\right\}
\right)  =I\right\}  \subseteq\mu\left(  S_{k}\times S_{n-k}\right)  .
\]
In other words, the map $\mu$ is surjective.

So the map $\mu$ is both injective and surjective. In other words, $\mu$ is
bijective. In other words, $\mu$ is a bijection. In other words, the map%
\begin{align*}
S_{k}\times S_{n-k}  &  \rightarrow\left\{  \tau\in S_{n}\ \mid\ \tau\left(
\left\{  1,2,\ldots,k\right\}  \right)  =I\right\}  ,\\
\left(  \alpha,\beta\right)   &  \mapsto\sigma_{I,\alpha,\beta}%
\end{align*}
is a bijection\footnote{since $\mu$ is the map
\begin{align*}
S_{k}\times S_{n-k}  &  \rightarrow\left\{  \tau\in S_{n}\ \mid\ \tau\left(
\left\{  1,2,\ldots,k\right\}  \right)  =I\right\}  ,\\
\left(  \alpha,\beta\right)   &  \mapsto\sigma_{I,\alpha,\beta}%
\end{align*}
}. This solves Exercise \ref{exe.Ialbe} \textbf{(c)}.
\end{proof}

\subsection{Solution to Exercise \ref{exe.perm.transX}}

We prepare for solving Exercise \ref{exe.perm.transX} by proving a variety of
lemmas of varying triviality (all of which are simple). Our first lemma is a
simple property of the permutations $t_{i,j}$ defined in Definition
\ref{def.transposX}:

\begin{lemma}
\label{lem.sol.perm.transX.tij1}Let $X$ be a set. Let $i$ and $j$ be two
distinct elements of $X$.

\textbf{(a)} We have $t_{i,j}\left(  i\right)  =j$.

\textbf{(b)} We have $t_{i,j}\left(  j\right)  =i$.

\textbf{(c)} We have $t_{i,j}\left(  k\right)  =k$ for each $k\in
X\setminus\left\{  i,j\right\}  $.

\textbf{(d)} We have $t_{i,j}\circ t_{i,j}=\operatorname*{id}$.
\end{lemma}

\begin{vershort}
\begin{proof}
[Proof of Lemma \ref{lem.sol.perm.transX.tij1}.]Lemma
\ref{lem.sol.perm.transX.tij1} follows from the definition of $t_{i,j}$.
\end{proof}
\end{vershort}

\begin{verlong}
\begin{proof}
[Proof of Lemma \ref{lem.sol.perm.transX.tij1}.]Definition \ref{def.transposX}
shows that $t_{i,j}$ is the permutation of $X$ which swaps $i$ with $j$ while
leaving all other elements of $X$ unchanged. In other words, the permutation
$t_{i,j}$ swaps $i$ with $j$, and leaves all elements of $X$ other than $i$
and $j$ unchanged.

The permutation $t_{i,j}$ swaps $i$ with $j$. In other words, it satisfies
$t_{i,j}\left(  i\right)  =j$ and $t_{i,j}\left(  j\right)  =i$. In
particular, we have $t_{i,j}\left(  i\right)  =j$. Hence, Lemma
\ref{lem.sol.perm.transX.tij1} \textbf{(a)} is proven. Also, $t_{i,j}\left(
j\right)  =i$. Thus, Lemma \ref{lem.sol.perm.transX.tij1} \textbf{(b)} is proven.

\textbf{(c)} Let $k\in X\setminus\left\{  i,j\right\}  $. We must prove that
$t_{i,j}\left(  k\right)  =k$.

We have $k\in X\setminus\left\{  i,j\right\}  $. In other words, $k$ is an
element of $X$ other than $i$ and $j$.

The permutation $t_{i,j}$ leaves all elements of $X$ other than $i$ and $j$
unchanged. In other words, it satisfies $t_{i,j}\left(  p\right)  =p$ whenever
$p$ is an element of $X$ other than $i$ and $j$. We can apply this to $p=k$
(since $k$ is an element of $X$ other than $i$ and $j$), and thus obtain
$t_{i,j}\left(  k\right)  =k$. Hence, Lemma \ref{lem.sol.perm.transX.tij1}
\textbf{(c)} is proven.

\textbf{(d)} Clearly, $t_{i,j}$ is a permutation of $X$, thus a bijection
$X\rightarrow X$. Hence, $t_{i,j}\circ t_{i,j}$ is a map $X\rightarrow X$ as well.

Let $k\in X$. We shall show that $\left(  t_{i,j}\circ t_{i,j}\right)  \left(
k\right)  =k$.

We have $k\in X$. We are in one of the following three cases:

\textit{Case 1:} We have $k=i$.

\textit{Case 2:} We have $k=j$.

\textit{Case 3:} We have neither $k=i$ nor $k=j$.

Let us first consider Case 1. In this case, we have $k=i$. Thus,%
\begin{align*}
\left(  t_{i,j}\circ t_{i,j}\right)  \left(  \underbrace{k}_{=i}\right)   &
=\left(  t_{i,j}\circ t_{i,j}\right)  \left(  i\right)  =t_{i,j}\left(
\underbrace{t_{i,j}\left(  i\right)  }_{\substack{=j\\\text{(by Lemma
\ref{lem.sol.perm.transX.tij1} \textbf{(a)})}}}\right)  =t_{i,j}\left(
j\right) \\
&  =i\ \ \ \ \ \ \ \ \ \ \left(  \text{by Lemma \ref{lem.sol.perm.transX.tij1}
\textbf{(b)}}\right) \\
&  =k.
\end{align*}
Thus, $\left(  t_{i,j}\circ t_{i,j}\right)  \left(  k\right)  =k$ is proven in
Case 1.

Let us now consider Case 2. In this case, we have $k=j$. Thus,%
\begin{align*}
\left(  t_{i,j}\circ t_{i,j}\right)  \left(  \underbrace{k}_{=j}\right)   &
=\left(  t_{i,j}\circ t_{i,j}\right)  \left(  j\right)  =t_{i,j}\left(
\underbrace{t_{i,j}\left(  j\right)  }_{\substack{=i\\\text{(by Lemma
\ref{lem.sol.perm.transX.tij1} \textbf{(b)})}}}\right)  =t_{i,j}\left(
i\right) \\
&  =j\ \ \ \ \ \ \ \ \ \ \left(  \text{by Lemma \ref{lem.sol.perm.transX.tij1}
\textbf{(a)}}\right) \\
&  =k.
\end{align*}
Thus, $\left(  t_{i,j}\circ t_{i,j}\right)  \left(  k\right)  =k$ is proven in
Case 2.

Let us finally consider Case 3. In this case, we have neither $k=i$ nor $k=j$.
Hence, we have $k\notin\left\{  i,j\right\}  $%
\ \ \ \ \footnote{\textit{Proof.} Assume the contrary. Thus, $k\in\left\{
i,j\right\}  $. In other words, $k$ is either $i$ or $j$. In other words, we
have either $k=i$ or $k=j$. This contradicts the fact that we have neither
$k=i$ nor $k=j$. This contradiction shows that our assumption was wrong.
Qed.}. Combining $k\in X$ with $k\notin\left\{  i,j\right\}  $, we obtain
$k\in X\setminus\left\{  i,j\right\}  $. Hence, $t_{i,j}\left(  k\right)  =k$
(by Lemma \ref{lem.sol.perm.transX.tij1} \textbf{(c)}). Now,
\[
\left(  t_{i,j}\circ t_{i,j}\right)  \left(  k\right)  =t_{i,j}\left(
\underbrace{t_{i,j}\left(  k\right)  }_{=k}\right)  =t_{i,j}\left(  k\right)
=k.
\]
Thus, $\left(  t_{i,j}\circ t_{i,j}\right)  \left(  k\right)  =k$ is proven in
Case 3.

We now have proven $\left(  t_{i,j}\circ t_{i,j}\right)  \left(  k\right)  =k$
in each of the three Cases 1, 2 and 3. Since these three Cases cover all
possibilities, we thus conclude that $\left(  t_{i,j}\circ t_{i,j}\right)
\left(  k\right)  =k$ always holds. Thus, $\left(  t_{i,j}\circ t_{i,j}%
\right)  \left(  k\right)  =k=\operatorname*{id}\left(  k\right)  $.

Now, forget that we fixed $k$. We thus have proven that $\left(  t_{i,j}\circ
t_{i,j}\right)  \left(  k\right)  =\operatorname*{id}\left(  k\right)  $ for
each $k\in X$. Hence, $t_{i,j}\circ t_{i,j}=\operatorname*{id}$ (since both
$t_{i,j}\circ t_{i,j}$ and $\operatorname*{id}$ are maps $X\rightarrow X$).
This proves Lemma \ref{lem.sol.perm.transX.tij1} \textbf{(d)}.
\end{proof}
\end{verlong}

\begin{lemma}
\label{lem.sol.perm.transX.tij-conj}Let $X$ and $Y$ be two sets. Let $i$ and
$j$ be two distinct elements of $X$. Let $f:X\rightarrow Y$ be a bijection.

\textbf{(a)} Then, $f\left(  i\right)  $ and $f\left(  j\right)  $ are two
distinct elements of $Y$, so that the transposition $t_{f\left(  i\right)
,f\left(  j\right)  }$ of $Y$ is well-defined.

\textbf{(b)} This transposition satisfies $t_{f\left(  i\right)  ,f\left(
j\right)  }=f\circ t_{i,j}\circ f^{-1}$.
\end{lemma}

\begin{vershort}
\begin{proof}
[Proof of Lemma \ref{lem.sol.perm.transX.tij-conj}.]Part \textbf{(a)} follows
from the injectivity of $f$. Part \textbf{(b)} is verified by directly
checking that $t_{f\left(  i\right)  ,f\left(  j\right)  }\left(  y\right)
=\left(  f\circ t_{i,j}\circ f^{-1}\right)  \left(  y\right)  $ for each $y\in
Y$.
\end{proof}
\end{vershort}

\begin{verlong}
\begin{proof}
[Proof of Lemma \ref{lem.sol.perm.transX.tij-conj}.]The map $f$ is bijective
(since $f$ is a bijection), and therefore injective. If we had $f\left(
i\right)  =f\left(  j\right)  $, then we would have $i=j$ (since $f$ is
injective), which would contradict the fact that $i$ and $j$ are distinct.
Thus, we cannot have $f\left(  i\right)  =f\left(  j\right)  $. Hence, we have
$f\left(  i\right)  \neq f\left(  j\right)  $. In other words, $f\left(
i\right)  $ and $f\left(  j\right)  $ are two distinct elements of $Y$. Hence,
the transposition $t_{f\left(  i\right)  ,f\left(  j\right)  }$ of $Y$ is
well-defined. This proves Lemma \ref{lem.sol.perm.transX.tij-conj}
\textbf{(a)}.

\textbf{(b)} Let $y\in Y$. We shall show that $t_{f\left(  i\right)  ,f\left(
j\right)  }\left(  y\right)  =\left(  f\circ t_{i,j}\circ f^{-1}\right)
\left(  y\right)  $.

We are in one of the following three cases:

\textit{Case 1:} We have $y=f\left(  i\right)  $.

\textit{Case 2:} We have $y=f\left(  j\right)  $.

\textit{Case 3:} We have neither $y=f\left(  i\right)  $ nor $y=f\left(
j\right)  $.

Let us first consider Case 1. In this case, we have $y=f\left(  i\right)  $.
Thus, $t_{f\left(  i\right)  ,f\left(  j\right)  }\left(  \underbrace{y}%
_{=f\left(  i\right)  }\right)  =t_{f\left(  i\right)  ,f\left(  j\right)
}\left(  f\left(  i\right)  \right)  =f\left(  j\right)  $ (by Lemma
\ref{lem.sol.perm.transX.tij1} \textbf{(a)} (applied to $Y$, $f\left(
i\right)  $ and $f\left(  j\right)  $ instead of $X$, $i$ and $j$)). But Lemma
\ref{lem.sol.perm.transX.tij1} \textbf{(a)} yields $t_{i,j}\left(  i\right)
=j$. Now,%
\begin{align*}
\left(  f\circ t_{i,j}\circ f^{-1}\right)  \left(  y\right)   &  =f\left(
t_{i,j}\left(  f^{-1}\left(  \underbrace{y}_{=f\left(  i\right)  }\right)
\right)  \right) \\
&  =f\left(  t_{i,j}\left(  \underbrace{f^{-1}\left(  f\left(  i\right)
\right)  }_{=i}\right)  \right)  =f\left(  \underbrace{t_{i,j}\left(
i\right)  }_{=j}\right)  =f\left(  j\right)  .
\end{align*}
Comparing this with $t_{f\left(  i\right)  ,f\left(  j\right)  }\left(
y\right)  =f\left(  j\right)  $, we obtain $t_{f\left(  i\right)  ,f\left(
j\right)  }\left(  y\right)  =\left(  f\circ t_{i,j}\circ f^{-1}\right)
\left(  y\right)  $. Thus, $t_{f\left(  i\right)  ,f\left(  j\right)  }\left(
y\right)  =\left(  f\circ t_{i,j}\circ f^{-1}\right)  \left(  y\right)  $ is
proven in Case 1.

Let us now consider Case 2. In this case, we have $y=f\left(  j\right)  $.
Thus, $t_{f\left(  i\right)  ,f\left(  j\right)  }\left(  \underbrace{y}%
_{=f\left(  j\right)  }\right)  =t_{f\left(  i\right)  ,f\left(  j\right)
}\left(  f\left(  j\right)  \right)  =f\left(  i\right)  $ (by Lemma
\ref{lem.sol.perm.transX.tij1} \textbf{(b)} (applied to $Y$, $f\left(
i\right)  $ and $f\left(  j\right)  $ instead of $X$, $i$ and $j$)). But Lemma
\ref{lem.sol.perm.transX.tij1} \textbf{(b)} yields $t_{i,j}\left(  j\right)
=i$. Now,%
\begin{align*}
\left(  f\circ t_{i,j}\circ f^{-1}\right)  \left(  y\right)   &  =f\left(
t_{i,j}\left(  f^{-1}\left(  \underbrace{y}_{=f\left(  j\right)  }\right)
\right)  \right) \\
&  =f\left(  t_{i,j}\left(  \underbrace{f^{-1}\left(  f\left(  j\right)
\right)  }_{=j}\right)  \right)  =f\left(  \underbrace{t_{i,j}\left(
j\right)  }_{=i}\right)  =f\left(  i\right)  .
\end{align*}
Comparing this with $t_{f\left(  i\right)  ,f\left(  j\right)  }\left(
y\right)  =f\left(  i\right)  $, we obtain $t_{f\left(  i\right)  ,f\left(
j\right)  }\left(  y\right)  =\left(  f\circ t_{i,j}\circ f^{-1}\right)
\left(  y\right)  $. Thus, $t_{f\left(  i\right)  ,f\left(  j\right)  }\left(
y\right)  =\left(  f\circ t_{i,j}\circ f^{-1}\right)  \left(  y\right)  $ is
proven in Case 2.

Let us finally consider Case 3. In this case, we have neither $y=f\left(
i\right)  $ nor $y=f\left(  j\right)  $. Hence, we have $y\notin\left\{
f\left(  i\right)  ,f\left(  j\right)  \right\}  $%
\ \ \ \ \footnote{\textit{Proof.} Assume the contrary. Thus, $y\in\left\{
f\left(  i\right)  ,f\left(  j\right)  \right\}  $. In other words, $y$ is
either $f\left(  i\right)  $ or $f\left(  j\right)  $. In other words, we have
either $y=f\left(  i\right)  $ or $y=f\left(  j\right)  $. This contradicts
the fact that we have neither $y=f\left(  i\right)  $ nor $y=f\left(
j\right)  $. This contradiction shows that our assumption was wrong. Qed.}.
Combining $y\in Y$ with $y\notin\left\{  f\left(  i\right)  ,f\left(
j\right)  \right\}  $, we obtain $y\in Y\setminus\left\{  f\left(  i\right)
,f\left(  j\right)  \right\}  $. Hence, Lemma \ref{lem.sol.perm.transX.tij1}
\textbf{(c)} (applied to $Y$, $f\left(  i\right)  $, $f\left(  j\right)  $ and
$y$ instead of $X$, $i$, $j$ and $k$) shows that $t_{f\left(  i\right)
,f\left(  j\right)  }\left(  y\right)  =y$.

On the other hand, $f^{-1}\left(  y\right)  \notin\left\{  i,j\right\}
$\ \ \ \ \footnote{\textit{Proof.} Assume the contrary. Thus, $f^{-1}\left(
y\right)  \in\left\{  i,j\right\}  $. But $f\left(  f^{-1}\left(  y\right)
\right)  =y$, so that $y=f\left(  \underbrace{f^{-1}\left(  y\right)  }%
_{\in\left\{  i,j\right\}  }\right)  \in f\left(  \left\{  i,j\right\}
\right)  =\left\{  f\left(  i\right)  ,f\left(  j\right)  \right\}  $. This
contradicts $y\notin\left\{  f\left(  i\right)  ,f\left(  j\right)  \right\}
$. This contradiction shows that our assumption was wrong. Qed.}. Combining
$f^{-1}\left(  y\right)  \in X$ with $f^{-1}\left(  y\right)  \notin\left\{
i,j\right\}  $, we obtain $f^{-1}\left(  y\right)  \in X\setminus\left\{
i,j\right\}  $. Thus, Lemma \ref{lem.sol.perm.transX.tij1} \textbf{(c)}
(applied to $k=f^{-1}\left(  y\right)  $) yields $t_{i,j}\left(  f^{-1}\left(
y\right)  \right)  =f^{-1}\left(  y\right)  $. Now,
\[
\left(  f\circ t_{i,j}\circ f^{-1}\right)  \left(  y\right)  =f\left(
\underbrace{t_{i,j}\left(  f^{-1}\left(  y\right)  \right)  }_{=f^{-1}\left(
y\right)  }\right)  =f\left(  f^{-1}\left(  y\right)  \right)  =y.
\]
Comparing this with $t_{f\left(  i\right)  ,f\left(  j\right)  }\left(
y\right)  =y$, we obtain $t_{f\left(  i\right)  ,f\left(  j\right)  }\left(
y\right)  =\left(  f\circ t_{i,j}\circ f^{-1}\right)  \left(  y\right)  $.
Thus, $t_{f\left(  i\right)  ,f\left(  j\right)  }\left(  y\right)  =\left(
f\circ t_{i,j}\circ f^{-1}\right)  \left(  y\right)  $ is proven in Case 3.

We now have proven $t_{f\left(  i\right)  ,f\left(  j\right)  }\left(
y\right)  =\left(  f\circ t_{i,j}\circ f^{-1}\right)  \left(  y\right)  $ in
each of the three Cases 1, 2 and 3. Since these three Cases cover all
possibilities, we thus conclude that $t_{f\left(  i\right)  ,f\left(
j\right)  }\left(  y\right)  =\left(  f\circ t_{i,j}\circ f^{-1}\right)
\left(  y\right)  $ always holds.

Now, forget that we fixed $y$. We thus have proven that $t_{f\left(  i\right)
,f\left(  j\right)  }\left(  y\right)  =\left(  f\circ t_{i,j}\circ
f^{-1}\right)  \left(  y\right)  $ for each $y\in Y$. Hence, $t_{f\left(
i\right)  ,f\left(  j\right)  }=f\circ t_{i,j}\circ f^{-1}$ (since both
$t_{f\left(  i\right)  ,f\left(  j\right)  }$ and $f\circ t_{i,j}\circ f^{-1}$
are maps $Y\rightarrow Y$). Thus, Lemma \ref{lem.sol.perm.transX.tij-conj}
\textbf{(b)} is proven.
\end{proof}
\end{verlong}

\begin{corollary}
\label{cor.sol.perm.transX.tij-sign}Let $X$ be a finite set. Let $i$ and $j$
be two distinct elements of $X$. Then, $\left(  -1\right)  ^{t_{i,j}}=-1$.
\end{corollary}

\begin{proof}
[Proof of Corollary \ref{cor.sol.perm.transX.tij-sign}.]Define an
$n\in\mathbb{N}$ by $n=\left\vert X\right\vert $. (This is well-defined, since
$X$ is a finite set.) Then, $\left\vert X\right\vert =n=\left\vert \left\{
1,2,\ldots,n\right\}  \right\vert $ (since $\left\vert \left\{  1,2,\ldots
,n\right\}  \right\vert =n$). Thus, there exists a bijection $\phi
:X\rightarrow\left\{  1,2,\ldots,n\right\}  $. Consider such a $\phi$.

Let $\sigma$ be the permutation $t_{i,j}$ of $X$. Thus, $\sigma=t_{i,j}$.
Consider the number $\left(  -1\right)  _{\phi}^{\sigma}$ defined as in
Exercise \ref{exe.ps4.2}. Then, the definition of $\left(  -1\right)
^{\sigma}$ (in Exercise \ref{exe.ps4.2}) yields $\left(  -1\right)  ^{\sigma
}=\left(  -1\right)  _{\phi}^{\sigma}=\left(  -1\right)  ^{\phi\circ
\sigma\circ\phi^{-1}}$ (by the definition of $\left(  -1\right)  _{\phi
}^{\sigma}$).

Lemma \ref{lem.sol.perm.transX.tij-conj} \textbf{(a)} (applied to $Y=\left\{
1,2,\ldots,n\right\}  $ and $f=\phi$) shows that $\phi\left(  i\right)  $ and
$\phi\left(  j\right)  $ are two distinct elements of $\left\{  1,2,\ldots
,n\right\}  $, so that the transposition $t_{\phi\left(  i\right)
,\phi\left(  j\right)  }$ of $\left\{  1,2,\ldots,n\right\}  $ is
well-defined. Lemma \ref{lem.sol.perm.transX.tij-conj} \textbf{(b)} (applied
to $Y=\left\{  1,2,\ldots,n\right\}  $ and $f=\phi$) shows that this
transposition satisfies $t_{\phi\left(  i\right)  ,\phi\left(  j\right)
}=\phi\circ t_{i,j}\circ\phi^{-1}$. But Exercise \ref{exe.ps4.1ab}
\textbf{(b)} (applied to $\phi\left(  i\right)  $ and $\phi\left(  j\right)  $
instead of $i$ and $j$) shows that $\left(  -1\right)  ^{t_{\phi\left(
i\right)  ,\phi\left(  j\right)  }}=-1$. In view of $t_{\phi\left(  i\right)
,\phi\left(  j\right)  }=\phi\circ t_{i,j}\circ\phi^{-1}$, this rewrites as
$\left(  -1\right)  ^{\phi\circ t_{i,j}\circ\phi^{-1}}=-1$.

Now,%
\begin{align*}
\left(  -1\right)  ^{\sigma}  &  =\left(  -1\right)  _{\phi}^{\sigma}=\left(
-1\right)  ^{\phi\circ\sigma\circ\phi^{-1}}=\left(  -1\right)  ^{\phi\circ
t_{i,j}\circ\phi^{-1}}\ \ \ \ \ \ \ \ \ \ \left(  \text{since }\sigma
=t_{i,j}\right) \\
&  =-1.
\end{align*}
In view of $\sigma=t_{i,j}$, this rewrites as $\left(  -1\right)  ^{t_{i,j}%
}=-1$. This proves Corollary \ref{cor.sol.perm.transX.tij-sign}.
\end{proof}

\begin{lemma}
\label{lem.sol.perm.transX.conj-prod}Let $X$ and $Y$ be two sets. Let
$f:X\rightarrow Y$ be a bijection. Let $p\in\mathbb{N}$. Let $\alpha
_{1},\alpha_{2},\ldots,\alpha_{p}$ be $p$ maps $Y\rightarrow Y$. For each
$i\in\left\{  1,2,\ldots,p\right\}  $, define a map $\beta_{i}:X\rightarrow X$
by $\beta_{i}=f^{-1}\circ\alpha_{i}\circ f$. Then,%
\[
\beta_{1}\circ\beta_{2}\circ\cdots\circ\beta_{p}=f^{-1}\circ\left(  \alpha
_{1}\circ\alpha_{2}\circ\cdots\circ\alpha_{p}\right)  \circ f.
\]

\end{lemma}

\begin{vershort}
\begin{proof}
[Proof of Lemma \ref{lem.sol.perm.transX.conj-prod}.]We claim that%
\begin{equation}
\beta_{1}\circ\beta_{2}\circ\cdots\circ\beta_{m}=f^{-1}\circ\left(  \alpha
_{1}\circ\alpha_{2}\circ\cdots\circ\alpha_{m}\right)  \circ f
\label{pf.lem.sol.perm.transX.conj-prod.short.goal}%
\end{equation}
for each $m\in\left\{  0,1,\ldots,p\right\}  $.

Indeed, (\ref{pf.lem.sol.perm.transX.conj-prod.short.goal}) can be proven by a
straightforward induction on $m$ (using the definition of $\beta_{i}$).

Now, applying (\ref{pf.lem.sol.perm.transX.conj-prod.short.goal}) to $m=p$, we
find $\beta_{1}\circ\beta_{2}\circ\cdots\circ\beta_{p}=f^{-1}\circ\left(
\alpha_{1}\circ\alpha_{2}\circ\cdots\circ\alpha_{p}\right)  \circ f$. This
proves Lemma \ref{lem.sol.perm.transX.conj-prod}.
\end{proof}
\end{vershort}

\begin{verlong}
\begin{proof}
[Proof of Lemma \ref{lem.sol.perm.transX.conj-prod}.]We claim that%
\begin{equation}
\beta_{1}\circ\beta_{2}\circ\cdots\circ\beta_{m}=f^{-1}\circ\left(  \alpha
_{1}\circ\alpha_{2}\circ\cdots\circ\alpha_{m}\right)  \circ f
\label{pf.lem.sol.perm.transX.conj-prod.goal}%
\end{equation}
for each $m\in\left\{  0,1,\ldots,p\right\}  $.

[\textit{Proof of (\ref{pf.lem.sol.perm.transX.conj-prod.goal}):} We shall
prove (\ref{pf.lem.sol.perm.transX.conj-prod.goal}) by induction on $m$:

\textit{Induction base:} We have $\beta_{1}\circ\beta_{2}\circ\cdots\circ
\beta_{0}=\left(  \text{empty composition of maps }X\rightarrow X\right)
=\operatorname*{id}\nolimits_{X}$ and $\alpha_{1}\circ\alpha_{2}\circ
\cdots\circ\alpha_{0}=\left(  \text{empty composition of maps }Y\rightarrow
Y\right)  =\operatorname*{id}\nolimits_{Y}$. Hence,%
\[
f^{-1}\circ\underbrace{\left(  \alpha_{1}\circ\alpha_{2}\circ\cdots\circ
\alpha_{0}\right)  }_{=\operatorname*{id}\nolimits_{Y}}\circ f=f^{-1}\circ
f=\operatorname*{id}\nolimits_{X}.
\]
Comparing this with $\beta_{1}\circ\beta_{2}\circ\cdots\circ\beta
_{0}=\operatorname*{id}\nolimits_{X}$, we obtain $\beta_{1}\circ\beta_{2}%
\circ\cdots\circ\beta_{0}=f^{-1}\circ\left(  \alpha_{1}\circ\alpha_{2}%
\circ\cdots\circ\alpha_{0}\right)  \circ f$. In other words,
(\ref{pf.lem.sol.perm.transX.conj-prod.goal}) holds for $m=0$. This completes
the induction base.

\textit{Induction step:} Let $M\in\left\{  0,1,\ldots,p\right\}  $ be
positive. Assume that (\ref{pf.lem.sol.perm.transX.conj-prod.goal}) holds for
$m=M-1$. We must prove that (\ref{pf.lem.sol.perm.transX.conj-prod.goal})
holds for $m=M$.

We have assumed that (\ref{pf.lem.sol.perm.transX.conj-prod.goal}) holds for
$m=M-1$. In other words, we have%
\[
\beta_{1}\circ\beta_{2}\circ\cdots\circ\beta_{M-1}=f^{-1}\circ\left(
\alpha_{1}\circ\alpha_{2}\circ\cdots\circ\alpha_{M-1}\right)  \circ f.
\]

We have $M\neq0$ (since $M$ is positive). Combining $M\in\left\{
0,1,\ldots,p\right\}  $ with $M\neq0$, we find $M\in\left\{  0,1,\ldots
,p\right\}  \setminus\left\{  0\right\}  =\left\{  1,2,\ldots,p\right\}  $.
Thus, the definition of $\beta_{M}$ yields $\beta_{M}=f^{-1}\circ\alpha
_{M}\circ f$. Now, $M$ is positive, so that%
\begin{align*}
\beta_{1}\circ\beta_{2}\circ\cdots\circ\beta_{M}  &  =\underbrace{\left(
\beta_{1}\circ\beta_{2}\circ\cdots\circ\beta_{M-1}\right)  }_{=f^{-1}%
\circ\left(  \alpha_{1}\circ\alpha_{2}\circ\cdots\circ\alpha_{M-1}\right)
\circ f}\circ\underbrace{\beta_{M}}_{=f^{-1}\circ\alpha_{M}\circ f}\\
&  =f^{-1}\circ\left(  \alpha_{1}\circ\alpha_{2}\circ\cdots\circ\alpha
_{M-1}\right)  \circ\underbrace{f\circ f^{-1}}_{=\operatorname*{id}%
\nolimits_{Y}}\circ\alpha_{M}\circ f\\
&  =f^{-1}\circ\underbrace{\left(  \alpha_{1}\circ\alpha_{2}\circ\cdots
\circ\alpha_{M-1}\right)  \circ\alpha_{M}}_{=\alpha_{1}\circ\alpha_{2}%
\circ\cdots\circ\alpha_{M}}\circ f\\
&  =f^{-1}\circ\left(  \alpha_{1}\circ\alpha_{2}\circ\cdots\circ\alpha
_{M}\right)  \circ f.
\end{align*}
In other words, (\ref{pf.lem.sol.perm.transX.conj-prod.goal}) holds for $m=M$.
This completes the induction step. Thus,
(\ref{pf.lem.sol.perm.transX.conj-prod.goal}) is proven by induction.]

Now, $p\in\left\{  0,1,\ldots,p\right\}  $ (since $p\in\mathbb{N}$). Hence,
(\ref{pf.lem.sol.perm.transX.conj-prod.goal}) (applied to $m=p$) yields
$\beta_{1}\circ\beta_{2}\circ\cdots\circ\beta_{p}=f^{-1}\circ\left(
\alpha_{1}\circ\alpha_{2}\circ\cdots\circ\alpha_{p}\right)  \circ f$. This
proves Lemma \ref{lem.sol.perm.transX.conj-prod}.
\end{proof}
\end{verlong}

The following fact is clear from the definitions:

\begin{proposition}
\label{prop.sol.perm.transX.sitra}Let $n\in\mathbb{N}$. Let $i\in\left\{
1,2,\ldots,n-1\right\}  $. Then, $s_{i}=t_{i,i+1}$.
\end{proposition}

\begin{verlong}
\begin{proof}
[Proof of Proposition \ref{prop.sol.perm.transX.sitra}.]We have $i\in\left\{
1,2,\ldots,n-1\right\}  $, so that $i+1\in\left\{  2,3,\ldots,n\right\}
\subseteq\left\{  1,2,\ldots,n\right\}  $. Also, $i\in\left\{  1,2,\ldots
,n-1\right\}  \subseteq\left\{  1,2,\ldots,n\right\}  $. Thus, $i$ and $i+1$
are two elements of $\left\{  1,2,\ldots,n\right\}  $. These two elements $i$
and $i+1$ are distinct (since $i\neq i+1$).

Recall that $t_{i,i+1}$ is defined to be the permutation in $S_{n}$ which
swaps $i$ with $i+1$ while leaving all other elements of $\left\{
1,2,\ldots,n\right\}  $ unchanged. In other words,%
\begin{align}
t_{i,i+1}  &  =\left(  \text{the permutation in }S_{n}\text{ which swaps
}i\text{ with }i+1\text{ while leaving}\right. \nonumber\\
&  \ \ \ \ \ \ \ \ \ \ \left.  \text{all other elements of }\left\{
1,2,\ldots,n\right\}  \text{ unchanged}\right)  .
\label{pf.prop.sol.perm.transX.sitra.tij}%
\end{align}

On the other hand, $s_{i}$ is defined as the permutation in $S_{n}$ that swaps
$i$ with $i+1$ but leaves all other numbers unchanged. In other words, $s_{i}$
is defined as the permutation in $S_{n}$ that swaps $i$ with $i+1$ but leaves
all other elements of $\left\{  1,2,\ldots,n\right\}  $ unchanged. In other
words,
\begin{align*}
s_{i}  &  =\left(  \text{the permutation in }S_{n}\text{ which swaps }i\text{
with }i+1\text{ while leaving}\right. \\
&  \ \ \ \ \ \ \ \ \ \ \left.  \text{all other elements of }\left\{
1,2,\ldots,n\right\}  \text{ unchanged}\right)  .
\end{align*}
Comparing this with (\ref{pf.prop.sol.perm.transX.sitra.tij}), we obtain
$s_{i}=t_{i,i+1}$. This proves Proposition \ref{prop.sol.perm.transX.sitra}.
\end{proof}
\end{verlong}

\begin{proposition}
\label{prop.sol.perm.transX.product}Let $X$ be a finite set. Let $\tau$ be any
permutation of $X$. Then, $\tau$ can be written as a composition of finitely
many transpositions of $X$.
\end{proposition}

\begin{vershort}
\begin{proof}
[Proof of Proposition \ref{prop.sol.perm.transX.product}.]Define an
$n\in\mathbb{N}$ by $n=\left\vert X\right\vert $. (This is well-defined, since
$X$ is a finite set.) Then, there exists a bijection $f:X\rightarrow\left\{
1,2,\ldots,n\right\}  $. Consider such a $f$.

Let $Y$ be the set $\left\{  1,2,\ldots,n\right\}  $. Thus, $Y=\left\{
1,2,\ldots,n\right\}  $. Hence, $f$ is a bijection $X\rightarrow Y$ (since $f$
is a bijection $X\rightarrow\left\{  1,2,\ldots,n\right\}  $).

But $\tau$ is a permutation of $X$. In other words, $\tau$ is a bijection
$X\rightarrow X$.

The maps $f$, $\tau$ and $f^{-1}$ are bijections. Hence, their composition
$f\circ\tau\circ f^{-1}:\left\{  1,2,\ldots,n\right\}  \rightarrow\left\{
1,2,\ldots,n\right\}  $ is also a bijection. Denote this bijection $f\circ
\tau\circ f^{-1}$ by $\sigma$. Thus, $\sigma=f\circ\tau\circ f^{-1}$.

We know that $\sigma$ is a bijection $\left\{  1,2,\ldots,n\right\}
\rightarrow\left\{  1,2,\ldots,n\right\}  $. In other words, $\sigma$ is a
permutation of the set $\left\{  1,2,\ldots,n\right\}  $. In other words,
$\sigma\in S_{n}$. Hence, Exercise \ref{exe.ps2.2.4} \textbf{(b)} shows that
$\sigma$ can be written as a composition of several permutations of the form
$s_{k}$ (with $k\in\left\{  1,2,\ldots,n-1\right\}  $). In other words, there
exist some $p\in\mathbb{N}$ and some elements $k_{1},k_{2},\ldots,k_{p}$ of
$\left\{  1,2,\ldots,n-1\right\}  $ such that $\sigma=s_{k_{1}}\circ s_{k_{2}%
}\circ\cdots\circ s_{k_{p}}$. Consider this $p$ and these $k_{1},k_{2}%
,\ldots,k_{p}$.

For each $i\in\left\{  1,2,\ldots,p\right\}  $, the map $s_{k_{i}}$ is a
bijection $Y\rightarrow Y$\ \ \ \ \footnote{\textit{Proof.} Let $i\in\left\{
1,2,\ldots,p\right\}  $. Then, $s_{k_{i}}\in S_{n}$. In other words,
$s_{k_{i}}$ is a permutation of the set $\left\{  1,2,\ldots,n\right\}  $. In
other words, $s_{k_{i}}$ is a permutation of the set $Y$ (since $Y=\left\{
1,2,\ldots,n\right\}  $). In other words, $s_{k_{i}}$ is a bijection
$Y\rightarrow Y$. Qed.}. Hence, for each $i\in\left\{  1,2,\ldots,p\right\}
$, we can define a map $\beta_{i}:X\rightarrow X$ by $\beta_{i}=f^{-1}\circ
s_{k_{i}}\circ f$ (because $f$ is a bijection $X\rightarrow Y$). Consider
these $\beta_{i}$. Lemma \ref{lem.sol.perm.transX.conj-prod} (applied to
$\alpha_{i}=s_{k_{i}}$) yields that%
\[
\beta_{1}\circ\beta_{2}\circ\cdots\circ\beta_{p}=f^{-1}\circ
\underbrace{\left(  s_{k_{1}}\circ s_{k_{2}}\circ\cdots\circ s_{k_{p}}\right)
}_{=\sigma=f\circ\tau\circ f^{-1}}\circ f=\underbrace{f^{-1}\circ
f}_{=\operatorname*{id}}\circ\tau\circ\underbrace{f^{-1}\circ f}%
_{=\operatorname*{id}}=\tau.
\]

But for each $i\in\left\{  1,2,\ldots,p\right\}  $, the map $\beta_{i}$ is a
transposition of $X$\ \ \ \ \footnote{\textit{Proof.} Let $i\in\left\{
1,2,\ldots,p\right\}  $. Thus, the definition of $\beta_{i}$ yields $\beta
_{i}=f^{-1}\circ s_{k_{i}}\circ f$.
\par
We have $k_{i}\in\left\{  1,2,\ldots,n-1\right\}  $. Hence, Proposition
\ref{prop.sol.perm.transX.sitra} (applied to $k_{i}$ instead of $i$) shows
that $s_{k_{i}}=t_{k_{i},k_{i}+1}$.
\par
But $k_{i}\in\left\{  1,2,\ldots,n-1\right\}  $. Hence, $k_{i}$ and $k_{i}+1$
are two elements of $\left\{  1,2,\ldots,n\right\}  =Y$. Thus, $f^{-1}\left(
k_{i}\right)  $ and $f^{-1}\left(  k_{i}+1\right)  $ are two elements of $X$
(since $f:X\rightarrow Y$ is a bijection). These two elements $f^{-1}\left(
k_{i}\right)  $ and $f^{-1}\left(  k_{i}+1\right)  $ are distinct (because if
we had $f^{-1}\left(  k_{i}\right)  =f^{-1}\left(  k_{i}+1\right)  $, then we
would have $k_{i}=k_{i}+1$, which is absurd). Thus, the map $t_{f^{-1}\left(
k_{i}\right)  ,f^{-1}\left(  k_{i}+1\right)  }$ is a transposition of $X$ (by
the definition of a \textquotedblleft transposition of $X$\textquotedblright).
\par
But Lemma \ref{lem.sol.perm.transX.tij-conj} \textbf{(b)} (applied to
$f^{-1}\left(  k_{i}\right)  $ and $f^{-1}\left(  k_{i}+1\right)  $ instead of
$i$ and $j$) shows that the transposition $t_{f\left(  f^{-1}\left(
k_{i}\right)  \right)  ,f\left(  f^{-1}\left(  k_{i}+1\right)  \right)  }$
satisfies $t_{f\left(  f^{-1}\left(  k_{i}\right)  \right)  ,f\left(
f^{-1}\left(  k_{i}+1\right)  \right)  }=f\circ t_{f^{-1}\left(  k_{i}\right)
,f^{-1}\left(  k_{i}+1\right)  }\circ f^{-1}$. Thus,%
\begin{align*}
f\circ t_{f^{-1}\left(  k_{i}\right)  ,f^{-1}\left(  k_{i}+1\right)  }\circ
f^{-1}  &  =t_{f\left(  f^{-1}\left(  k_{i}\right)  \right)  ,f\left(
f^{-1}\left(  k_{i}+1\right)  \right)  }\\
&  =t_{k_{i},k_{i}+1}\ \ \ \ \ \ \ \ \ \ \left(  \text{since }f\left(
f^{-1}\left(  k_{i}\right)  \right)  =k_{i}\text{ and }f\left(  f^{-1}\left(
k_{i}+1\right)  \right)  =k_{i}+1\right) \\
&  =s_{k_{i}}\ \ \ \ \ \ \ \ \ \ \left(  \text{since }s_{k_{i}}=t_{k_{i}%
,k_{i}+1}\right)  .
\end{align*}
Hence,%
\[
f^{-1}\circ\underbrace{f\circ t_{f^{-1}\left(  k_{i}\right)  ,f^{-1}\left(
k_{i}+1\right)  }\circ f^{-1}}_{=s_{k_{i}}}\circ f=f^{-1}\circ s_{k_{i}}\circ
f=\beta_{i}\ \ \ \ \ \ \ \ \ \ \left(  \text{since }\beta_{i}=f^{-1}\circ
s_{k_{i}}\circ f\right)  .
\]
Thus, $\beta_{i}=\underbrace{f^{-1}\circ f}_{=\operatorname*{id}\nolimits_{X}%
}\circ t_{f^{-1}\left(  k_{i}\right)  ,f^{-1}\left(  k_{i}+1\right)  }%
\circ\underbrace{f^{-1}\circ f}_{=\operatorname*{id}\nolimits_{X}}%
=t_{f^{-1}\left(  k_{i}\right)  ,f^{-1}\left(  k_{i}+1\right)  }$.
\par
But we know that the map $t_{f^{-1}\left(  k_{i}\right)  ,f^{-1}\left(
k_{i}+1\right)  }$ is a transposition of $X$. In view of $\beta_{i}%
=t_{f^{-1}\left(  k_{i}\right)  ,f^{-1}\left(  k_{i}+1\right)  }$, this
rewrites as follows: The map $\beta_{i}$ is a transposition of $X$. Qed.}. In
other words, $\beta_{1},\beta_{2},\ldots,\beta_{p}$ are transpositions of $X$.
Thus, the permutation $\beta_{1}\circ\beta_{2}\circ\cdots\circ\beta_{p}$ is a
composition of finitely many transpositions of $X$. In view of $\beta_{1}%
\circ\beta_{2}\circ\cdots\circ\beta_{p}=\tau$, this rewrites as follows: The
permutation $\tau$ is a composition of finitely many transpositions of $X$.
This proves Proposition \ref{prop.sol.perm.transX.product}.
\end{proof}
\end{vershort}

\begin{verlong}
\begin{proof}
[Proof of Proposition \ref{prop.sol.perm.transX.product}.]Define an
$n\in\mathbb{N}$ by $n=\left\vert X\right\vert $. (This is well-defined, since
$X$ is a finite set.) Then, $\left\vert X\right\vert =n=\left\vert \left\{
1,2,\ldots,n\right\}  \right\vert $ (since $\left\vert \left\{  1,2,\ldots
,n\right\}  \right\vert =n$). Thus, there exists a bijection $f:X\rightarrow
\left\{  1,2,\ldots,n\right\}  $. Consider such a $f$.

Let $Y$ be the set $\left\{  1,2,\ldots,n\right\}  $. Thus, $Y=\left\{
1,2,\ldots,n\right\}  $. Hence, $f$ is a bijection $X\rightarrow Y$ (since $f$
is a bijection $X\rightarrow\left\{  1,2,\ldots,n\right\}  $).

The map $f$ is a bijection $X\rightarrow\left\{  1,2,\ldots,n\right\}  $.
Thus, the inverse $f^{-1}$ of $f$ is well-defined and is a bijection $\left\{
1,2,\ldots,n\right\}  \rightarrow X$.

The map $\tau$ is a permutation of $X$. In other words, the map $\tau$ is a
bijection $X\rightarrow X$.

The maps $f$, $\tau$ and $f^{-1}$ are bijections. Hence, their composition
$f\circ\tau\circ f^{-1}:\left\{  1,2,\ldots,n\right\}  \rightarrow\left\{
1,2,\ldots,n\right\}  $ is also a bijection. Denote this bijection $f\circ
\tau\circ f^{-1}$ by $\sigma$. Thus, $\sigma=f\circ\tau\circ f^{-1}$.

We know that $\sigma$ is a bijection $\left\{  1,2,\ldots,n\right\}
\rightarrow\left\{  1,2,\ldots,n\right\}  $. In other words, $\sigma$ is a
permutation of the set $\left\{  1,2,\ldots,n\right\}  $. In other words,
$\sigma\in S_{n}$ (since $S_{n}$ is the set of all permutations of the set
$\left\{  1,2,\ldots,n\right\}  $). Hence, Exercise \ref{exe.ps2.2.4}
\textbf{(b)} shows that $\sigma$ can be written as a composition of several
permutations of the form $s_{k}$ (with $k\in\left\{  1,2,\ldots,n-1\right\}
$). In other words, there exist some $p\in\mathbb{N}$ and some elements
$k_{1},k_{2},\ldots,k_{p}$ of $\left\{  1,2,\ldots,n-1\right\}  $ such that
$\sigma=s_{k_{1}}\circ s_{k_{2}}\circ\cdots\circ s_{k_{p}}$. Consider this $p$
and these $k_{1},k_{2},\ldots,k_{p}$.

For each $i\in\left\{  1,2,\ldots,p\right\}  $, the map $s_{k_{i}}$ is a
bijection $Y\rightarrow Y$\ \ \ \ \footnote{\textit{Proof.} Let $i\in\left\{
1,2,\ldots,p\right\}  $. Then, $s_{k_{i}}\in S_{n}$. In other words,
$s_{k_{i}}$ is a permutation of the set $\left\{  1,2,\ldots,n\right\}  $
(since $S_{n}$ is the set of all permutations of the set $\left\{
1,2,\ldots,n\right\}  $). In other words, $s_{k_{i}}$ is a permutation of the
set $Y$ (since $Y=\left\{  1,2,\ldots,n\right\}  $). In other words,
$s_{k_{i}}$ is a bijection $Y\rightarrow Y$. Qed.}. Hence, for each
$i\in\left\{  1,2,\ldots,p\right\}  $, we can define a map $\beta
_{i}:X\rightarrow X$ by $\beta_{i}=f^{-1}\circ s_{k_{i}}\circ f$ (because $f$
is a bijection $X\rightarrow Y$). Consider these $\beta_{i}$. Lemma
\ref{lem.sol.perm.transX.conj-prod} (applied to $\alpha_{i}=s_{k_{i}}$) yields
that%
\[
\beta_{1}\circ\beta_{2}\circ\cdots\circ\beta_{p}=f^{-1}\circ
\underbrace{\left(  s_{k_{1}}\circ s_{k_{2}}\circ\cdots\circ s_{k_{p}}\right)
}_{=\sigma=f\circ\tau\circ f^{-1}}\circ f=\underbrace{f^{-1}\circ
f}_{=\operatorname*{id}}\circ\tau\circ\underbrace{f^{-1}\circ f}%
_{=\operatorname*{id}}=\tau.
\]
In other words, $\tau=\beta_{1}\circ\beta_{2}\circ\cdots\circ\beta_{p}$.

But for each $i\in\left\{  1,2,\ldots,p\right\}  $, the map $\beta_{i}$ is a
transposition of $X$\ \ \ \ \footnote{\textit{Proof.} Let $i\in\left\{
1,2,\ldots,p\right\}  $. Thus, the definition of $\beta_{i}$ yields $\beta
_{i}=f^{-1}\circ s_{k_{i}}\circ f$.
\par
We have $k_{i}\in\left\{  1,2,\ldots,n-1\right\}  $. Hence, Proposition
\ref{prop.sol.perm.transX.sitra} (applied to $k_{i}$ instead of $i$) shows
that $s_{k_{i}}=t_{k_{i},k_{i}+1}$.
\par
But $k_{i}\in\left\{  1,2,\ldots,n-1\right\}  \subseteq\left\{  1,2,\ldots
,n\right\}  =Y$. Also, from $k_{i}\in\left\{  1,2,\ldots,n-1\right\}  $, we
obtain $k_{i}+1\in\left\{  2,3,\ldots,n\right\}  \subseteq\left\{
1,2,\ldots,n\right\}  =Y$. Thus, $k_{i}$ and $k_{i}+1$ are two elements of
$Y$.
\par
Recall that $f$ is a bijection $X\rightarrow Y$. Thus, $f^{-1}$ is a
well-defined map $Y\rightarrow X$. Hence, the elements $f^{-1}\left(
k_{i}\right)  $ and $f^{-1}\left(  k_{i}+1\right)  $ of $X$ are defined (since
$k_{i}$ and $k_{i}+1$ are two elements of $Y$).
\par
Assume (for the sake of contradiction) that $f^{-1}\left(  k_{i}\right)
=f^{-1}\left(  k_{i}+1\right)  $. Thus, $f\left(  \underbrace{f^{-1}\left(
k_{i}\right)  }_{=f^{-1}\left(  k_{i}+1\right)  }\right)  =f\left(
f^{-1}\left(  k_{i}+1\right)  \right)  =k_{i}+1$. Comparing this with
$f\left(  f^{-1}\left(  k_{i}\right)  \right)  =k_{i}$, we obtain $k_{i}%
=k_{i}+1$. This contradicts $k_{i}\neq k_{i}+1$. This contradiction shows that
our assumption (that $f^{-1}\left(  k_{i}\right)  =f^{-1}\left(
k_{i}+1\right)  $) was wrong. Hence, we don't have $f^{-1}\left(
k_{i}\right)  =f^{-1}\left(  k_{i}+1\right)  $. In other words, we have
$f^{-1}\left(  k_{i}\right)  \neq f^{-1}\left(  k_{i}+1\right)  $. In other
words, the two elements $f^{-1}\left(  k_{i}\right)  $ and $f^{-1}\left(
k_{i}+1\right)  $ of $X$ are distinct. Thus, the map $t_{f^{-1}\left(
k_{i}\right)  ,f^{-1}\left(  k_{i}+1\right)  }$ is a transposition of $X$ (by
the definition of a \textquotedblleft transposition of $X$\textquotedblright).
\par
But Lemma \ref{lem.sol.perm.transX.tij-conj} \textbf{(b)} (applied to
$f^{-1}\left(  k_{i}\right)  $ and $f^{-1}\left(  k_{i}+1\right)  $ instead of
$i$ and $j$) shows that the transposition $t_{f\left(  f^{-1}\left(
k_{i}\right)  \right)  ,f\left(  f^{-1}\left(  k_{i}+1\right)  \right)  }$
satisfies $t_{f\left(  f^{-1}\left(  k_{i}\right)  \right)  ,f\left(
f^{-1}\left(  k_{i}+1\right)  \right)  }=f\circ t_{f^{-1}\left(  k_{i}\right)
,f^{-1}\left(  k_{i}+1\right)  }\circ f^{-1}$. Thus,%
\begin{align*}
f\circ t_{f^{-1}\left(  k_{i}\right)  ,f^{-1}\left(  k_{i}+1\right)  }\circ
f^{-1}  &  =t_{f\left(  f^{-1}\left(  k_{i}\right)  \right)  ,f\left(
f^{-1}\left(  k_{i}+1\right)  \right)  }\\
&  =t_{k_{i},k_{i}+1}\ \ \ \ \ \ \ \ \ \ \left(  \text{since }f\left(
f^{-1}\left(  k_{i}\right)  \right)  =k_{i}\text{ and }f\left(  f^{-1}\left(
k_{i}+1\right)  \right)  =k_{i}+1\right) \\
&  =s_{k_{i}}\ \ \ \ \ \ \ \ \ \ \left(  \text{since }s_{k_{i}}=t_{k_{i}%
,k_{i}+1}\right)  .
\end{align*}
Hence,%
\[
f^{-1}\circ\underbrace{f\circ t_{f^{-1}\left(  k_{i}\right)  ,f^{-1}\left(
k_{i}+1\right)  }\circ f^{-1}}_{=s_{k_{i}}}\circ f=f^{-1}\circ s_{k_{i}}\circ
f=\beta_{i}\ \ \ \ \ \ \ \ \ \ \left(  \text{since }\beta_{i}=f^{-1}\circ
s_{k_{i}}\circ f\right)  .
\]
Thus,%
\[
\beta_{i}=\underbrace{f^{-1}\circ f}_{=\operatorname*{id}\nolimits_{X}}\circ
t_{f^{-1}\left(  k_{i}\right)  ,f^{-1}\left(  k_{i}+1\right)  }\circ
\underbrace{f^{-1}\circ f}_{=\operatorname*{id}\nolimits_{X}}=t_{f^{-1}\left(
k_{i}\right)  ,f^{-1}\left(  k_{i}+1\right)  }.
\]
\par
But we know that the map $t_{f^{-1}\left(  k_{i}\right)  ,f^{-1}\left(
k_{i}+1\right)  }$ is a transposition of $X$. In view of $\beta_{i}%
=t_{f^{-1}\left(  k_{i}\right)  ,f^{-1}\left(  k_{i}+1\right)  }$, this
rewrites as follows: The map $\beta_{i}$ is a transposition of $X$. Qed.}. In
other words, $\beta_{1},\beta_{2},\ldots,\beta_{p}$ are transpositions of $X$.
Thus, the permutation $\beta_{1}\circ\beta_{2}\circ\cdots\circ\beta_{p}$ is a
composition of finitely many transpositions of $X$. In view of $\tau=\beta
_{1}\circ\beta_{2}\circ\cdots\circ\beta_{p}$, this rewrites as follows: The
permutation $\tau$ is a composition of finitely many transpositions of $X$.
This proves Proposition \ref{prop.sol.perm.transX.product}.
\end{proof}
\end{verlong}

\begin{proposition}
\label{prop.sol.perm.transX.sign.prod-of-many}Let $X$ be a finite set. Let
$k\in\mathbb{N}$. Let $\sigma_{1},\sigma_{2},\ldots,\sigma_{k}$ be $k$
permutations of $X$. Then,%
\[
\left(  -1\right)  ^{\sigma_{1}\circ\sigma_{2}\circ\cdots\circ\sigma_{k}%
}=\left(  -1\right)  ^{\sigma_{1}}\cdot\left(  -1\right)  ^{\sigma_{2}}%
\cdot\cdots\cdot\left(  -1\right)  ^{\sigma_{k}}.
\]

\end{proposition}

\begin{vershort}
\begin{proof}
[Proof of Proposition \ref{prop.sol.perm.transX.sign.prod-of-many}.]This is
analogous to the proof of Proposition \ref{prop.sign.prod-of-many}, using
Exercise \ref{exe.ps4.2} \textbf{(c)}.
\end{proof}
\end{vershort}

\begin{verlong}
\begin{proof}
[Proof of Proposition \ref{prop.sol.perm.transX.sign.prod-of-many}.]We claim
that%
\begin{equation}
\left(  -1\right)  ^{\sigma_{1}\circ\sigma_{2}\circ\cdots\circ\sigma_{m}%
}=\left(  -1\right)  ^{\sigma_{1}}\cdot\left(  -1\right)  ^{\sigma_{2}}%
\cdot\cdots\cdot\left(  -1\right)  ^{\sigma_{m}}
\label{pf.prop.sol.perm.transX.sign.prod-of-many.goal}%
\end{equation}
for each $m\in\left\{  0,1,\ldots,k\right\}  $.

[\textit{Proof of (\ref{pf.prop.sol.perm.transX.sign.prod-of-many.goal}):}] We
will prove (\ref{pf.prop.sol.perm.transX.sign.prod-of-many.goal}) by induction
over $m$:

\textit{Induction base:} We have $\sigma_{1}\circ\sigma_{2}\circ\cdots
\circ\sigma_{0}=\left(  \text{composition of }0\text{ permutations of
}X\right)  =\operatorname*{id}$. Hence, $\left(  -1\right)  ^{\sigma_{1}%
\circ\sigma_{2}\circ\cdots\circ\sigma_{0}}=\left(  -1\right)
^{\operatorname*{id}}=1$ (by Exercise \ref{exe.ps4.2} \textbf{(b)}). On the
other hand, from $\left(  -1\right)  ^{\sigma_{1}}\cdot\left(  -1\right)
^{\sigma_{2}}\cdot\cdots\cdot\left(  -1\right)  ^{\sigma_{0}}=\left(
\text{product of }0\text{ integers}\right)  =1$. Compared with $\left(
-1\right)  ^{\sigma_{1}\circ\sigma_{2}\circ\cdots\circ\sigma_{0}}=1$, this
yields $\left(  -1\right)  ^{\sigma_{1}\circ\sigma_{2}\circ\cdots\circ
\sigma_{0}}=\left(  -1\right)  ^{\sigma_{1}}\cdot\left(  -1\right)
^{\sigma_{2}}\cdot\cdots\cdot\left(  -1\right)  ^{\sigma_{0}}$. In other
words, (\ref{pf.prop.sol.perm.transX.sign.prod-of-many.goal}) holds for $m=0$.
The induction base is thus complete.

\textit{Induction step:} Let $M\in\mathbb{N}$. Assume that
(\ref{pf.prop.sol.perm.transX.sign.prod-of-many.goal}) is proven for $m=M$. We
need to prove (\ref{pf.prop.sol.perm.transX.sign.prod-of-many.goal}) for
$m=M+1$.

We have assumed that (\ref{pf.prop.sol.perm.transX.sign.prod-of-many.goal}) is
proven for $m=M$. In other words, we have $\left(  -1\right)  ^{\sigma
_{1}\circ\sigma_{2}\circ\cdots\circ\sigma_{M}}=\left(  -1\right)  ^{\sigma
_{1}}\cdot\left(  -1\right)  ^{\sigma_{2}}\cdot\cdots\cdot\left(  -1\right)
^{\sigma_{M}}$.

But $\sigma_{1}\circ\sigma_{2}\circ\cdots\circ\sigma_{M+1}=\left(  \sigma
_{1}\circ\sigma_{2}\circ\cdots\circ\sigma_{M}\right)  \circ\sigma_{M+1}$, so
that%
\begin{align*}
\left(  -1\right)  ^{\sigma_{1}\circ\sigma_{2}\circ\cdots\circ\sigma_{M+1}}
&  =\left(  -1\right)  ^{\left(  \sigma_{1}\circ\sigma_{2}\circ\cdots
\circ\sigma_{M}\right)  \circ\sigma_{M+1}}=\underbrace{\left(  -1\right)
^{\sigma_{1}\circ\sigma_{2}\circ\cdots\circ\sigma_{M}}}_{=\left(  -1\right)
^{\sigma_{1}}\cdot\left(  -1\right)  ^{\sigma_{2}}\cdot\cdots\cdot\left(
-1\right)  ^{\sigma_{M}}}\cdot\left(  -1\right)  ^{\sigma_{M+1}}\\
&  \ \ \ \ \ \ \ \ \ \ \left(  \text{by Exercise \ref{exe.ps4.2} \textbf{(c)}
(applied to }\sigma=\sigma_{1}\circ\sigma_{2}\circ\cdots\circ\sigma_{M}\text{
and }\tau=\sigma_{M+1}\text{)}\right) \\
&  =\left(  \left(  -1\right)  ^{\sigma_{1}}\cdot\left(  -1\right)
^{\sigma_{2}}\cdot\cdots\cdot\left(  -1\right)  ^{\sigma_{M}}\right)
\cdot\left(  -1\right)  ^{\sigma_{M+1}}\\
&  =\left(  -1\right)  ^{\sigma_{1}}\cdot\left(  -1\right)  ^{\sigma_{2}}%
\cdot\cdots\cdot\left(  -1\right)  ^{\sigma_{M+1}}.
\end{align*}
In other words, (\ref{pf.prop.sol.perm.transX.sign.prod-of-many.goal}) holds
for $m=M+1$. This completes the induction step, and so the induction proof of
(\ref{pf.prop.sol.perm.transX.sign.prod-of-many.goal}) is complete.]

Now, $k\in\left\{  0,1,\ldots,k\right\}  $ (since $k\in\mathbb{N}$). Hence,
(\ref{pf.prop.sol.perm.transX.sign.prod-of-many.goal}) (applied to $m=k$)
yields $\left(  -1\right)  ^{\sigma_{1}\circ\sigma_{2}\circ\cdots\circ
\sigma_{k}}=\left(  -1\right)  ^{\sigma_{1}}\cdot\left(  -1\right)
^{\sigma_{2}}\cdot\cdots\cdot\left(  -1\right)  ^{\sigma_{k}}$. Thus,
Proposition \ref{prop.sol.perm.transX.sign.prod-of-many} is proven.
\end{proof}
\end{verlong}

\begin{corollary}
\label{cor.sol.perm.transX.sign.prod-of-trans}Let $X$ be a finite set. Let
$k\in\mathbb{N}$. Let $\sigma$ be a permutation of $X$.\ Assume that $\sigma$
can be written as a composition of $k$ transpositions of $X$. Then, $\left(
-1\right)  ^{\sigma}=\left(  -1\right)  ^{k}$.
\end{corollary}

\begin{proof}
[Proof of Corollary \ref{cor.sol.perm.transX.sign.prod-of-trans}.]We know that
$\sigma$ can be written as a composition of $k$ transpositions of $X$. In
other words, there exist $k$ transpositions $u_{1},u_{2},\ldots,u_{k}$ of $X$
such that $\sigma=u_{1}\circ u_{2}\circ\cdots\circ u_{k}$. Consider these
$u_{1},u_{2},\ldots,u_{k}$.

\begin{vershort}
Hence, $u_{1},u_{2},\ldots,u_{k}$ are $k$ permutations of $X$.
\end{vershort}

\begin{verlong}
Now, $u_{1},u_{2},\ldots,u_{k}$ are $k$ transpositions of $X$, thus $k$
permutations of $X$ (since each transposition of $X$ is a permutation of $X$).
\end{verlong}

Thus, Proposition \ref{prop.sol.perm.transX.sign.prod-of-many} (applied to
$\sigma_{i}=u_{i}$) yields
\[
\left(  -1\right)  ^{u_{1}\circ u_{2}\circ\cdots\circ u_{k}}=\left(
-1\right)  ^{u_{1}}\cdot\left(  -1\right)  ^{u_{2}}\cdot\cdots\cdot\left(
-1\right)  ^{u_{k}}=\prod_{p=1}^{k}\left(  -1\right)  ^{u_{p}}.
\]

Each $p\in\left\{  1,2,\ldots,k\right\}  $ satisfies $\left(  -1\right)
^{u_{p}}=-1$\ \ \ \ \footnote{\textit{Proof.} Let $p\in\left\{  1,2,\ldots
,k\right\}  $. Then, $u_{p}$ is a transposition of $X$ (since $u_{1}%
,u_{2},\ldots,u_{k}$ are $k$ transpositions of $X$). In other words, $u_{p}$
has the form $u_{p}=t_{i,j}$ where $i$ and $j$ are two distinct elements of
$X$ (by the definition of a \textquotedblleft transposition of $X$%
\textquotedblright). Consider these $i$ and $j$. From $u_{p}=t_{i,j}$, we
obtain $\left(  -1\right)  ^{u_{p}}=\left(  -1\right)  ^{t_{i,j}}=-1$ (by
Corollary \ref{cor.sol.perm.transX.tij-sign}). Qed.}. Thus, $\prod_{p=1}%
^{k}\underbrace{\left(  -1\right)  ^{u_{p}}}_{=-1}=\prod_{p=1}^{k}\left(
-1\right)  =\left(  -1\right)  ^{k}$.

From $\sigma=u_{1}\circ u_{2}\circ\cdots\circ u_{k}$, we obtain $\left(
-1\right)  ^{\sigma}=\left(  -1\right)  ^{u_{1}\circ u_{2}\circ\cdots\circ
u_{k}}=\prod_{p=1}^{k}\left(  -1\right)  ^{u_{p}}=\left(  -1\right)  ^{k}$.
This proves Corollary \ref{cor.sol.perm.transX.sign.prod-of-trans}.
\end{proof}

\begin{vershort}
\begin{proof}
[Solution to Exercise \ref{exe.perm.transX}.]Exercise \ref{exe.perm.transX}
\textbf{(a)} is precisely Corollary \ref{cor.sol.perm.transX.tij-sign}.
Exercise \ref{exe.perm.transX} \textbf{(b)} is precisely Proposition
\ref{prop.sol.perm.transX.product}. Exercise \ref{exe.perm.transX}
\textbf{(c)} is precisely Corollary
\ref{cor.sol.perm.transX.sign.prod-of-trans}. Thus, Exercise
\ref{exe.perm.transX} is solved.
\end{proof}
\end{vershort}

\begin{verlong}
\begin{proof}
[Solution to Exercise \ref{exe.perm.transX}.]\textbf{(a)} Let $i$ and $j$ be
two distinct elements of $X$. Then, Corollary
\ref{cor.sol.perm.transX.tij-sign} yields $\left(  -1\right)  ^{t_{i,j}}=-1$.
This solves Exercise \ref{exe.perm.transX} \textbf{(a)}.

\textbf{(b)} If $\tau$ is any permutation of $X$, then $\tau$ can be written
as a composition of finitely many transpositions of $X$ (by Proposition
\ref{prop.sol.perm.transX.product}). In other words, any permutation $\tau$ of
$X$ can be written as a composition of finitely many transpositions of $X$.
This solves Exercise \ref{exe.perm.transX} \textbf{(b)}.

\textbf{(c)} Corollary \ref{cor.sol.perm.transX.sign.prod-of-trans} yields
$\left(  -1\right)  ^{\sigma}=\left(  -1\right)  ^{k}$. This solves Exercise
\ref{exe.perm.transX} \textbf{(c)}.
\end{proof}
\end{verlong}

\subsection{\label{sect.sol.perm.c=ttt}Solution to Exercise
\ref{exe.perm.c=ttt}}

Throughout this section, we shall use Definition \ref{def.transpos}. We shall
also use the notation $\left[  n\right]  $ for the set $\left\{
1,2,\ldots,n\right\}  $ whenever $n\in\mathbb{N}$.

\subsubsection{The \textquotedblleft moving lemmas\textquotedblright}

Before we start solving Exercise \ref{exe.perm.c=ttt}, let us state three
lemmas that help us understand how a composition of several maps acts on an
element. All three of these lemmas are intuitively obvious; the length of
their formal proofs is mainly due to the amount of bookkeeping they require.

\begin{lemma}
\label{lem.sol.sigmacrosstau.compose0}Let $X$ be a set. Let $m\in\mathbb{N}$.
Let $f_{1},f_{2},\ldots,f_{m}$ be $m$ maps from $X$ to $X$. Let $x\in X$.

Assume that%
\begin{equation}
f_{j}\left(  x\right)  =x\ \ \ \ \ \ \ \ \ \ \text{for each }j\in\left\{
1,2,\ldots,m\right\}  . \label{eq.lem.sol.sigmacrosstau.compose0.ass}%
\end{equation}
Then, $\left(  f_{m}\circ f_{m-1}\circ\cdots\circ f_{1}\right)  \left(
x\right)  =x$.
\end{lemma}

\begin{vershort}
\begin{proof}
[Proof of Lemma \ref{lem.sol.sigmacrosstau.compose0}.]Here is an informal
proof: Imagine the element $x$ undergoing the maps $f_{1},f_{2},\ldots,f_{m}$
in this order; the result is, of course, $\left(  f_{m}\circ f_{m-1}%
\circ\cdots\circ f_{1}\right)  \left(  x\right)  $. But let us look closer at
the step-by-step procedure. The element is initially $x$. Then, the maps
$f_{1},f_{2},\ldots,f_{m}$ are being applied to it in this order. The element
never changes in the process (because of
(\ref{eq.lem.sol.sigmacrosstau.compose0.ass})). Thus, the final result is
still $x$. This shows that $\left(  f_{m}\circ f_{m-1}\circ\cdots\circ
f_{1}\right)  \left(  x\right)  =x$.

If you wish, you can easily turn this argument into a rigorous proof: We claim
that each $k\in\left\{  0,1,\ldots,m\right\}  $ satisfies%
\begin{equation}
\left(  f_{k}\circ f_{k-1}\circ\cdots\circ f_{1}\right)  \left(  x\right)  =x.
\label{pf.lem.sol.sigmacrosstau.compose0.short.goal}%
\end{equation}
The proof of (\ref{pf.lem.sol.sigmacrosstau.compose0.short.goal}) is a
straightforward induction on $k$ (where
(\ref{eq.lem.sol.sigmacrosstau.compose0.ass}) is used in the induction step).
After (\ref{pf.lem.sol.sigmacrosstau.compose0.short.goal}) is proven, we can
apply (\ref{pf.lem.sol.sigmacrosstau.compose0.short.goal}) to $k=m$, and
conclude that $\left(  f_{m}\circ f_{m-1}\circ\cdots\circ f_{1}\right)
\left(  x\right)  =x$. Thus, Lemma \ref{lem.sol.sigmacrosstau.compose0} is proven.
\end{proof}
\end{vershort}

\begin{verlong}
\begin{proof}
[Proof of Lemma \ref{lem.sol.sigmacrosstau.compose0}.]We claim that each
$k\in\left\{  0,1,\ldots,m\right\}  $ satisfies%
\begin{equation}
\left(  f_{k}\circ f_{k-1}\circ\cdots\circ f_{1}\right)  \left(  x\right)  =x.
\label{pf.lem.sol.sigmacrosstau.compose0.goal}%
\end{equation}

[\textit{Proof of (\ref{pf.lem.sol.sigmacrosstau.compose0.goal}):} We shall
prove (\ref{pf.lem.sol.sigmacrosstau.compose0.goal}) by induction over $k$:

\textit{Induction base:} We have%
\[
\underbrace{\left(  f_{0}\circ f_{0-1}\circ\cdots\circ f_{1}\right)
}_{=\left(  \text{empty composition of maps }X\rightarrow X\right)
=\operatorname*{id}\nolimits_{X}}\left(  x\right)  =\operatorname*{id}%
\nolimits_{X}\left(  x\right)  =x.
\]
In other words, (\ref{pf.lem.sol.sigmacrosstau.compose0.goal}) holds for
$k=0$. This completes the induction base.

\textit{Induction step:} Let $K\in\left\{  0,1,\ldots,m\right\}  $ be
positive. Assume that (\ref{pf.lem.sol.sigmacrosstau.compose0.goal}) holds for
$k=K-1$. We must prove that (\ref{pf.lem.sol.sigmacrosstau.compose0.goal})
holds for $k=K$.

We have $K\neq0$ (since $K$ is positive). Combining this with $K\in\left\{
0,1,\ldots,m\right\}  $, we obtain $K\in\left\{  0,1,\ldots,m\right\}
\setminus\left\{  1\right\}  =\left\{  1,2,\ldots,m\right\}  $. Hence,
$K-1\in\left\{  0,1,\ldots,m-1\right\}  \subseteq\left\{  0,1,\ldots
,m\right\}  $. But we have assumed that
(\ref{pf.lem.sol.sigmacrosstau.compose0.goal}) holds for $k=K-1$. In other
words, we have $\left(  f_{K-1}\circ f_{\left(  K-1\right)  -1}\circ
\cdots\circ f_{1}\right)  \left(  x\right)  =x$.

But (\ref{eq.lem.sol.sigmacrosstau.compose0.ass}) (applied to $j=K$) yields
$f_{K}\left(  x\right)  =x$. Now,%
\begin{align*}
f_{K}\circ f_{K-1}\circ\cdots\circ f_{1}  &  =f_{K}\circ\underbrace{\left(
f_{K-1}\circ f_{K-2}\circ\cdots\circ f_{1}\right)  }_{=f_{K-1}\circ f_{\left(
K-1\right)  -1}\circ\cdots\circ f_{1}}\\
&  =f_{K}\circ\left(  f_{K-1}\circ f_{\left(  K-1\right)  -1}\circ\cdots\circ
f_{1}\right)  .
\end{align*}
Applying both sides of this equality to $x$, we find%
\begin{align*}
\left(  f_{K}\circ f_{K-1}\circ\cdots\circ f_{1}\right)  \left(  x\right)   &
=\left(  f_{K}\circ\left(  f_{K-1}\circ f_{\left(  K-1\right)  -1}\circ
\cdots\circ f_{1}\right)  \right)  \left(  x\right) \\
&  =f_{K}\left(  \underbrace{\left(  f_{K-1}\circ f_{\left(  K-1\right)
-1}\circ\cdots\circ f_{1}\right)  \left(  x\right)  }_{=x}\right)
=f_{K}\left(  x\right)  =x.
\end{align*}
In other words, (\ref{pf.lem.sol.sigmacrosstau.compose0.goal}) holds for
$k=K$. This completes the induction step. Thus,
(\ref{pf.lem.sol.sigmacrosstau.compose0.goal}) is proven.]

Now, $m\in\mathbb{N}$, so that $m\in\left\{  0,1,\ldots,m\right\}  $. Hence,
(\ref{pf.lem.sol.sigmacrosstau.compose0.goal}) (applied to $k=m$) yields
$\left(  f_{m}\circ f_{m-1}\circ\cdots\circ f_{1}\right)  \left(  x\right)
=x$. This proves Lemma \ref{lem.sol.sigmacrosstau.compose0}.
\end{proof}
\end{verlong}

\begin{lemma}
\label{lem.sol.sigmacrosstau.compose1}Let $X$ be a set. Let $m\in\mathbb{N}$.
Let $f_{1},f_{2},\ldots,f_{m}$ be $m$ maps from $X$ to $X$. Let $x$ and $y$ be
two elements of $X$.

Let $i\in\left\{  1,2,\ldots,m\right\}  $. Assume that $f_{i}\left(  x\right)
=y$. Assume further that%
\begin{equation}
f_{j}\left(  x\right)  =x\ \ \ \ \ \ \ \ \ \ \text{for each }j\in\left\{
1,2,\ldots,m\right\}  \text{ satisfying }j<i.
\label{eq.lem.sol.sigmacrosstau.compose1.ass1}%
\end{equation}
Assume also that%
\begin{equation}
f_{j}\left(  y\right)  =y\ \ \ \ \ \ \ \ \ \ \text{for each }j\in\left\{
1,2,\ldots,m\right\}  \text{ satisfying }j>i.
\label{eq.lem.sol.sigmacrosstau.compose1.ass2}%
\end{equation}
Then, $\left(  f_{m}\circ f_{m-1}\circ\cdots\circ f_{1}\right)  \left(
x\right)  =y$.
\end{lemma}

\begin{vershort}
\begin{proof}
[Proof of Lemma \ref{lem.sol.sigmacrosstau.compose1}.]Lemma
\ref{lem.sol.sigmacrosstau.compose1} is a more complicated variant of Lemma
\ref{lem.sol.sigmacrosstau.compose0}. Again, there is a quick informal proof:
Imagine the element $x$ undergoing the maps $f_{1},f_{2},\ldots,f_{m}$ in this
order; the result is, of course, $\left(  f_{m}\circ f_{m-1}\circ\cdots\circ
f_{1}\right)  \left(  x\right)  $. But let us look closer at the step-by-step
procedure. The element is initially $x$. Then, the maps $f_{1},f_{2}%
,\ldots,f_{m}$ are being applied to it in this order. Up until the map $f_{i}$
is applied, the element does not change (because of
(\ref{eq.lem.sol.sigmacrosstau.compose1.ass1})). Then, the map $f_{i}$ is
applied, and the element becomes $y$ (since $f_{i}\left(  x\right)  =y$). From
then on, the maps $f_{i+1},f_{i+2},\ldots,f_{m}$ again leave the element
unchanged (due to (\ref{eq.lem.sol.sigmacrosstau.compose1.ass2})). Thus, the
final result is $y$. This shows that $\left(  f_{m}\circ f_{m-1}\circ
\cdots\circ f_{1}\right)  \left(  x\right)  =y$.

Again, it is not hard to formalize this argument into a rigorous proof: For
each $p\in\mathbb{Z}$, we define an element $z_{p}\in X$ by%
\[
z_{p}=%
\begin{cases}
x, & \text{if }p<i;\\
y, & \text{if }p\geq i
\end{cases}
.
\]
These elements have the following property:%
\begin{equation}
f_{j}\left(  z_{j-1}\right)  =z_{j}\ \ \ \ \ \ \ \ \ \ \text{for each }%
j\in\left\{  1,2,\ldots,m\right\}  .
\label{pf.lem.sol.sigmacrosstau.compose1.short.rec}%
\end{equation}
(Indeed, this is proven mechanically by distinguishing the cases $j<i$, $j=i$
and $j>i$. In the case when $j<i$, the equality
(\ref{pf.lem.sol.sigmacrosstau.compose1.short.rec}) follows from
(\ref{eq.lem.sol.sigmacrosstau.compose1.ass1}). In the case when $j=i$, the
equality (\ref{pf.lem.sol.sigmacrosstau.compose1.short.rec}) follows from
$f_{i}\left(  x\right)  =y$. In the case when $j>i$, the equality
(\ref{pf.lem.sol.sigmacrosstau.compose1.short.rec}) follows from
(\ref{eq.lem.sol.sigmacrosstau.compose1.ass2}).)

Notice that $i\in\left\{  1,2,\ldots,m\right\}  $, so that $0<i$ and $m\geq
i$. The definition of $z_{0}$ yields $z_{0}=x$ (since $0<i$), whereas the
definition of $z_{m}$ yields $z_{m}=y$ (since $m\geq i$).

Now, we claim that each $k\in\left\{  0,1,\ldots,m\right\}  $ satisfies%
\begin{equation}
\left(  f_{k}\circ f_{k-1}\circ\cdots\circ f_{1}\right)  \left(  x\right)
=z_{k}. \label{pf.lem.sol.sigmacrosstau.compose1.short.goal}%
\end{equation}
The proof of (\ref{pf.lem.sol.sigmacrosstau.compose1.short.goal}) is a
straightforward induction on $k$ (where $z_{0}=x$ is used in the induction
base, and (\ref{pf.lem.sol.sigmacrosstau.compose1.short.rec}) is used in the
induction step). After (\ref{pf.lem.sol.sigmacrosstau.compose1.short.goal}) is
proven, we can apply (\ref{pf.lem.sol.sigmacrosstau.compose1.short.goal}) to
$k=m$, and conclude that $\left(  f_{m}\circ f_{m-1}\circ\cdots\circ
f_{1}\right)  \left(  x\right)  =z_{m}=y$. Thus, Lemma
\ref{lem.sol.sigmacrosstau.compose1} is proven.
\end{proof}
\end{vershort}

\begin{verlong}
\begin{proof}
[Proof of Lemma \ref{lem.sol.sigmacrosstau.compose1}.]For each $p\in
\mathbb{Z}$, we define an element $z_{p}\in X$ by%
\[
z_{p}=%
\begin{cases}
x, & \text{if }p<i;\\
y, & \text{if }p\geq i
\end{cases}
.
\]
These elements have the following property:%
\begin{equation}
f_{j}\left(  z_{j-1}\right)  =z_{j}\ \ \ \ \ \ \ \ \ \ \text{for each }%
j\in\left\{  1,2,\ldots,m\right\}  .
\label{pf.lem.sol.sigmacrosstau.compose1.rec}%
\end{equation}

[\textit{Proof of (\ref{pf.lem.sol.sigmacrosstau.compose1.rec}):} Let
$j\in\left\{  1,2,\ldots,m\right\}  $. We want to prove
(\ref{pf.lem.sol.sigmacrosstau.compose1.rec}).

We are in one of the following three cases:

\textit{Case 1:} We have $j<i$.

\textit{Case 2:} We have $j=i$.

\textit{Case 3:} We have $j>i$.

Let us first consider Case 1. In this case, we have $j<i$, so that $j-1<j<i$.
The definition of $z_{j}$ yields $z_{j}=%
\begin{cases}
x, & \text{if }j<i;\\
y, & \text{if }j\geq i
\end{cases}
=x$ (since $j<i$). The definition of $z_{j-1}$ yields $z_{j-1}=%
\begin{cases}
x, & \text{if }j-1<i;\\
y, & \text{if }j-1\geq i
\end{cases}
=x$ (since $j-1<i$). Finally, (\ref{eq.lem.sol.sigmacrosstau.compose1.ass1})
shows that $f_{j}\left(  x\right)  =x$. Thus,%
\[
f_{j}\left(  \underbrace{z_{j-1}}_{=x}\right)  =f_{j}\left(  x\right)
=x=z_{j}.
\]
Hence, (\ref{pf.lem.sol.sigmacrosstau.compose1.rec}) is proven in Case 1.

Let us next consider Case 2. In this case, we have $j=i$, so that $j-1<j=i$
and $j=i\geq i$. The definition of $z_{j}$ yields $z_{j}=%
\begin{cases}
x, & \text{if }j<i;\\
y, & \text{if }j\geq i
\end{cases}
=y$ (since $j\geq i$). The definition of $z_{j-1}$ yields $z_{j-1}=%
\begin{cases}
x, & \text{if }j-1<i;\\
y, & \text{if }j-1\geq i
\end{cases}
=x$ (since $j-1<i$). Finally, $j=i$, so that $f_{j}\left(  x\right)
=f_{i}\left(  x\right)  =y$. Now,%
\[
f_{j}\left(  \underbrace{z_{j-1}}_{=x}\right)  =f_{j}\left(  x\right)
=y=z_{j}.
\]
Hence, (\ref{pf.lem.sol.sigmacrosstau.compose1.rec}) is proven in Case 2.

Let us finally consider Case 3. In this case, we have $j>i$, so that $j\geq
i+1$ (since $j$ and $i$ are integers). Hence, $j-1\geq i$, so that $j\geq
j-1\geq i$. The definition of $z_{j}$ yields $z_{j}=%
\begin{cases}
x, & \text{if }j<i;\\
y, & \text{if }j\geq i
\end{cases}
=y$ (since $j\geq i$). The definition of $z_{j-1}$ yields $z_{j-1}=%
\begin{cases}
x, & \text{if }j-1<i;\\
y, & \text{if }j-1\geq i
\end{cases}
=y$ (since $j-1\geq i$). Finally,
(\ref{eq.lem.sol.sigmacrosstau.compose1.ass2}) shows that $f_{j}\left(
y\right)  =y$. Thus,%
\[
f_{j}\left(  \underbrace{z_{j-1}}_{=y}\right)  =f_{j}\left(  y\right)
=y=z_{j}.
\]
Hence, (\ref{pf.lem.sol.sigmacrosstau.compose1.rec}) is proven in Case 3.

We have now proven (\ref{pf.lem.sol.sigmacrosstau.compose1.rec}) in each of
the three Cases 1, 2 and 3. Since these three Cases cover all possibilities,
we thus conclude that (\ref{pf.lem.sol.sigmacrosstau.compose1.rec}) always
holds. This completes the proof of
(\ref{pf.lem.sol.sigmacrosstau.compose1.rec}).]

We have $i\in\left\{  1,2,\ldots,m\right\}  $, so that $1\leq i\leq m$. Thus,
$m\geq i$. Now, we claim the following:

\begin{statement}
\textit{Observation 1:} We have $\left(  f_{g}\circ f_{g-1}\circ\cdots\circ
f_{1}\right)  \left(  x\right)  =z_{g}$ for each $g\in\left\{  0,1,\ldots
,m\right\}  $.
\end{statement}

[\textit{Proof of Observation 1:} We shall prove Observation 1 by induction on
$g$:

\textit{Induction base:} We have $i\geq1$ (since $i\in\left\{  1,2,\ldots
,m\right\}  $), so that $i\geq1>0$ and thus $0<i$. The definition of $z_{0}$
yields $z_{0}=%
\begin{cases}
x, & \text{if }0<i;\\
y, & \text{if }0\geq i
\end{cases}
=x$ (since $0<i$). Comparing this with%
\[
\underbrace{\left(  f_{0}\circ f_{0-1}\circ\cdots\circ f_{1}\right)
}_{\substack{=\left(  \text{empty composition of maps }X\rightarrow X\right)
\\=\operatorname*{id}\nolimits_{X}}}\left(  x\right)  =\operatorname*{id}%
\nolimits_{X}\left(  x\right)  =x,
\]
we obtain $\left(  f_{0}\circ f_{0-1}\circ\cdots\circ f_{1}\right)  \left(
x\right)  =z_{0}$. In other words, Observation 1 holds for $g=0$. This
completes the induction base.

\textit{Induction step:} Let $h\in\left\{  0,1,\ldots,m\right\}  $ be
positive. Assume that Observation 1 holds for $g=h-1$. We must then prove that
Observation 1 holds for $g=h$.

We have $h\neq0$ (since $h$ is positive). Combining this with $h\in\left\{
0,1,\ldots,m\right\}  $, we obtain $h\in\left\{  0,1,\ldots,m\right\}
\setminus\left\{  0\right\}  =\left\{  1,2,\ldots,m\right\}  $, so that
$h-1\in\left\{  0,1,\ldots,m-1\right\}  \subseteq\left\{  0,1,\ldots
,m\right\}  $.

But we assumed that Observation 1 holds for $g=h-1$. In other words, we have
$\left(  f_{h-1}\circ f_{\left(  h-1\right)  -1}\circ\cdots\circ f_{1}\right)
\left(  x\right)  =z_{h-1}$

Now, $h$ is positive. Thus,%
\begin{align*}
&  \underbrace{\left(  f_{h}\circ f_{h-1}\circ\cdots\circ f_{1}\right)
}_{=f_{h}\circ\left(  f_{h-1}\circ f_{h-2}\circ\cdots\circ f_{1}\right)
}\left(  x\right)  =\left(  f_{h}\circ\underbrace{\left(  f_{h-1}\circ
f_{h-2}\circ\cdots\circ f_{1}\right)  }_{=f_{h-1}\circ f_{\left(  h-1\right)
-1}\circ\cdots\circ f_{1}}\right)  \left(  x\right) \\
&  =\left(  f_{h}\circ\left(  f_{h-1}\circ f_{\left(  h-1\right)  -1}%
\circ\cdots\circ f_{1}\right)  \right)  \left(  x\right)  =f_{h}\left(
\underbrace{\left(  f_{h-1}\circ f_{\left(  h-1\right)  -1}\circ\cdots\circ
f_{1}\right)  \left(  x\right)  }_{=z_{h-1}}\right) \\
&  =f_{h}\left(  z_{h-1}\right)  =z_{h}\ \ \ \ \ \ \ \ \ \ \left(  \text{by
(\ref{pf.lem.sol.sigmacrosstau.compose1.rec}) (applied to }j=h\text{)}\right)
.
\end{align*}
In other words, Observation 1 holds for $g=h$. This completes the induction
step. Thus, Observation 1 is proven.]

We can now apply Observation 1 to $g=m$ (since $m\in\left\{  0,1,\ldots
,m\right\}  $ (because $m\in\mathbb{N}$)). We thus obtain%
\begin{align*}
\left(  f_{m}\circ f_{m-1}\circ\cdots\circ f_{1}\right)  \left(  x\right)   &
=z_{m}=%
\begin{cases}
x, & \text{if }m<i;\\
y, & \text{if }m\geq i
\end{cases}
\ \ \ \ \ \ \ \ \ \ \left(  \text{by the definition of }z_{m}\right) \\
&  =y\ \ \ \ \ \ \ \ \ \ \left(  \text{since }m\geq i\right)  .
\end{align*}
This proves Lemma \ref{lem.sol.sigmacrosstau.compose1}.
\end{proof}
\end{verlong}

The next lemma is a generalization of Lemma
\ref{lem.sol.sigmacrosstau.compose0} (and, with a bit more work, of Lemma
\ref{lem.sol.sigmacrosstau.compose1}):

\begin{lemma}
\label{lem.sol.sigmacrosstau.composeall}Let $X$ be a set. Let $m\in\mathbb{N}%
$. Let $f_{1},f_{2},\ldots,f_{m}$ be $m$ maps from $X$ to $X$. Let
$x_{0},x_{1},\ldots,x_{m}$ be $m+1$ elements of $X$.

Assume that%
\begin{equation}
f_{j}\left(  x_{j-1}\right)  =x_{j}\ \ \ \ \ \ \ \ \ \ \text{for each }%
j\in\left\{  1,2,\ldots,m\right\}  .
\label{eq.lem.sol.sigmacrosstau.composeall.ass}%
\end{equation}
Then, $\left(  f_{m}\circ f_{m-1}\circ\cdots\circ f_{1}\right)  \left(
x_{0}\right)  =x_{m}$.
\end{lemma}

\begin{vershort}
\begin{proof}
[Proof of Lemma \ref{lem.sol.sigmacrosstau.composeall}.]Here is an informal
proof: From (\ref{eq.lem.sol.sigmacrosstau.composeall.ass}), we obtain the $m$
equalities%
\[
f_{1}\left(  x_{0}\right)  =x_{1},\ \ \ \ \ \ \ \ \ \ f_{2}\left(
x_{1}\right)  =x_{2},\ \ \ \ \ \ \ \ \ \ ...,\ \ \ \ \ \ \ \ \ \ f_{m}\left(
x_{m-1}\right)  =x_{m}.
\]
Now, imagine the maps $f_{1},f_{2},\ldots,f_{m}$ being applied (in this order)
to the element $x_{0}$; the result is, of course, $\left(  f_{m}\circ
f_{m-1}\circ\cdots\circ f_{1}\right)  \left(  x_{0}\right)  $. But let us look
closer at this step-by-step procedure. We start with the element $x_{0}$.
Then, we apply the map $f_{1}$ to it, and it becomes $x_{1}$ (since
$f_{1}\left(  x_{0}\right)  =x_{1}$). Next, we apply the map $f_{2}$ to it
(i.e., to $x_{1}$, not to $x_{0}$), and it becomes $x_{2}$ (since
$f_{2}\left(  x_{1}\right)  =x_{2}$). Then, we apply the map $f_{3}$ to it,
and it becomes $x_{3}$ (for similar reasons). We then continue with
$f_{4},f_{5},\ldots,f_{m}$, obtaining the elements $x_{4},x_{5},\ldots,x_{m}$
in the process (for the same reasons); hence, the final result is $x_{m}$.
Thus, we know that the final result is simultaneously $\left(  f_{m}\circ
f_{m-1}\circ\cdots\circ f_{1}\right)  \left(  x_{0}\right)  $ and $x_{m}$.
Hence, $\left(  f_{m}\circ f_{m-1}\circ\cdots\circ f_{1}\right)  \left(
x_{0}\right)  =x_{m}$.

If you wish, you can easily turn this argument into a rigorous proof: We claim
that each $k\in\left\{  0,1,\ldots,m\right\}  $ satisfies%
\begin{equation}
\left(  f_{k}\circ f_{k-1}\circ\cdots\circ f_{1}\right)  \left(  x_{0}\right)
=x_{k}. \label{pf.lem.sol.sigmacrosstau.composeall.short.goal}%
\end{equation}
The proof of (\ref{pf.lem.sol.sigmacrosstau.composeall.short.goal}) is a
straightforward induction on $k$ (where
(\ref{eq.lem.sol.sigmacrosstau.composeall.ass}) is used in the induction
step). After (\ref{pf.lem.sol.sigmacrosstau.composeall.short.goal}) is proven,
we can apply (\ref{pf.lem.sol.sigmacrosstau.composeall.short.goal}) to $k=m$,
and conclude that $\left(  f_{m}\circ f_{m-1}\circ\cdots\circ f_{1}\right)
\left(  x_{0}\right)  =x_{m}$. Thus, Lemma
\ref{lem.sol.sigmacrosstau.composeall} is proven.
\end{proof}
\end{vershort}

\begin{verlong}
\begin{proof}
[Proof of Lemma \ref{lem.sol.sigmacrosstau.composeall}.]We claim that each
$k\in\left\{  0,1,\ldots,m\right\}  $ satisfies%
\begin{equation}
\left(  f_{k}\circ f_{k-1}\circ\cdots\circ f_{1}\right)  \left(  x_{0}\right)
=x_{k}. \label{pf.lem.sol.sigmacrosstau.composeall.goal}%
\end{equation}

[\textit{Proof of (\ref{pf.lem.sol.sigmacrosstau.composeall.goal}):} We shall
prove (\ref{pf.lem.sol.sigmacrosstau.composeall.goal}) by induction over $k$:

\textit{Induction base:} We have%
\[
\underbrace{\left(  f_{0}\circ f_{0-1}\circ\cdots\circ f_{1}\right)
}_{=\left(  \text{empty composition of maps }X\rightarrow X\right)
=\operatorname*{id}\nolimits_{X}}\left(  x_{0}\right)  =\operatorname*{id}%
\nolimits_{X}\left(  x_{0}\right)  =x_{0}.
\]
In other words, (\ref{pf.lem.sol.sigmacrosstau.composeall.goal}) holds for
$k=0$. This completes the induction base.

\textit{Induction step:} Let $K\in\left\{  0,1,\ldots,m\right\}  $ be
positive. Assume that (\ref{pf.lem.sol.sigmacrosstau.composeall.goal}) holds
for $k=K-1$. We must prove that
(\ref{pf.lem.sol.sigmacrosstau.composeall.goal}) holds for $k=K$.

We have $K\neq0$ (since $K$ is positive). Combining this with $K\in\left\{
0,1,\ldots,m\right\}  $, we obtain $K\in\left\{  0,1,\ldots,m\right\}
\setminus\left\{  1\right\}  =\left\{  1,2,\ldots,m\right\}  $. Hence,
$K-1\in\left\{  0,1,\ldots,m-1\right\}  \subseteq\left\{  0,1,\ldots
,m\right\}  $. But we have assumed that
(\ref{pf.lem.sol.sigmacrosstau.composeall.goal}) holds for $k=K-1$. In other
words, we have $\left(  f_{K-1}\circ f_{\left(  K-1\right)  -1}\circ
\cdots\circ f_{1}\right)  \left(  x_{0}\right)  =x_{K-1}$.

But (\ref{eq.lem.sol.sigmacrosstau.composeall.ass}) (applied to $j=K$) yields
$f_{K}\left(  x_{K-1}\right)  =x_{K}$. Now,%
\begin{align*}
f_{K}\circ f_{K-1}\circ\cdots\circ f_{1}  &  =f_{K}\circ\underbrace{\left(
f_{K-1}\circ f_{K-2}\circ\cdots\circ f_{1}\right)  }_{=f_{K-1}\circ f_{\left(
K-1\right)  -1}\circ\cdots\circ f_{1}}\\
&  =f_{K}\circ\left(  f_{K-1}\circ f_{\left(  K-1\right)  -1}\circ\cdots\circ
f_{1}\right)  .
\end{align*}
Applying both sides of this equality to $x_{0}$, we find%
\begin{align*}
\left(  f_{K}\circ f_{K-1}\circ\cdots\circ f_{1}\right)  \left(  x_{0}\right)
&  =\left(  f_{K}\circ\left(  f_{K-1}\circ f_{\left(  K-1\right)  -1}%
\circ\cdots\circ f_{1}\right)  \right)  \left(  x_{0}\right) \\
&  =f_{K}\left(  \underbrace{\left(  f_{K-1}\circ f_{\left(  K-1\right)
-1}\circ\cdots\circ f_{1}\right)  \left(  x_{0}\right)  }_{=x_{K-1}}\right)
=f_{K}\left(  x_{K-1}\right)  =x_{K}.
\end{align*}
In other words, (\ref{pf.lem.sol.sigmacrosstau.composeall.goal}) holds for
$k=K$. This completes the induction step. Thus,
(\ref{pf.lem.sol.sigmacrosstau.composeall.goal}) is proven.]

Now, $m\in\mathbb{N}$, so that $m\in\left\{  0,1,\ldots,m\right\}  $. Hence,
(\ref{pf.lem.sol.sigmacrosstau.composeall.goal}) (applied to $k=m$) yields
$\left(  f_{m}\circ f_{m-1}\circ\cdots\circ f_{1}\right)  \left(
x_{0}\right)  =x_{m}$. This proves Lemma
\ref{lem.sol.sigmacrosstau.composeall}.
\end{proof}
\end{verlong}

\subsubsection{Solving Exercise \ref{exe.perm.c=ttt}}

\begin{proof}
[Solution to Exercise \ref{exe.perm.c=ttt}.]We have $k\in\left\{
1,2,\ldots,n\right\}  $, thus $k\geq1$, hence $k-1\geq0$. Therefore,
$k-1\in\mathbb{N}$.

\begin{verlong}
For each $j\in\left\{  1,2,\ldots,k-1\right\}  $, the transposition
$t_{i_{j},i_{j+1}}$ in $S_{n}$ is well-defined\footnote{\textit{Proof.} Let
$j\in\left\{  1,2,\ldots,k-1\right\}  $. Thus, $j+1\in\left\{  2,3,\ldots
,k\right\}  \subseteq\left\{  1,2,\ldots,k\right\}  $ and $j\in\left\{
1,2,\ldots,k-1\right\}  \subseteq\left\{  1,2,\ldots,k\right\}  $. Hence, $j$
and $j+1$ are two elements of $\left\{  1,2,\ldots,k\right\}  $. Thus, $i_{j}$
and $i_{j+1}$ are two elements of $\left[  n\right]  $ (since $i_{1}%
,i_{2},\ldots,i_{k}$ are $k$ elements of $\left[  n\right]  $). Moreover,
$j\neq j+1$ and thus $i_{j}\neq i_{j+1}$ (since the $k$ elements $i_{1}%
,i_{2},\ldots,i_{k}$ are distinct). In other words, $i_{j}$ and $i_{j+1}$ are
distinct. Hence, $i_{j}$ and $i_{j+1}$ are two distinct elements of $\left[
n\right]  $. In other words, $i_{j}$ and $i_{j+1}$ are two distinct elements
of $\left\{  1,2,\ldots,n\right\}  $ (since $\left[  n\right]  =\left\{
1,2,\ldots,n\right\}  $). Thus, the transposition $t_{i_{j},i_{j+1}}$ in
$S_{n}$ is well-defined (according to Definition \ref{def.transpos}). Qed.}.
In other words, the $k-1$ transpositions $t_{i_{1},i_{2}},t_{i_{2},i_{3}%
},\ldots,t_{i_{k-1},i_{k}}$ in $S_{n}$ are well-defined.
\end{verlong}

We have defined $\operatorname*{cyc}\nolimits_{i_{1},i_{2},\ldots,i_{k}}$ to
be the permutation in $S_{n}$ which sends $i_{1},i_{2},\ldots,i_{k}$ to
$i_{2},i_{3},\ldots,i_{k},i_{1}$, respectively, while leaving all other
elements of $\left[  n\right]  $ fixed. Therefore:

\begin{itemize}
\item The permutation $\operatorname*{cyc}\nolimits_{i_{1},i_{2},\ldots,i_{k}%
}$ sends $i_{1},i_{2},\ldots,i_{k}$ to $i_{2},i_{3},\ldots,i_{k},i_{1}$,
respectively. In other words,%
\begin{equation}
\operatorname*{cyc}\nolimits_{i_{1},i_{2},\ldots,i_{k}}\left(  i_{p}\right)
=i_{p+1}\ \ \ \ \ \ \ \ \ \ \text{for every }p\in\left\{  1,2,\ldots
,k\right\}  , \label{sol.perm.c=ttt.cyc-manifest.1}%
\end{equation}
where $i_{k+1}$ means $i_{1}$.

\item The permutation $\operatorname*{cyc}\nolimits_{i_{1},i_{2},\ldots,i_{k}%
}$ leaves all other elements of $\left[  n\right]  $ fixed (where
\textquotedblleft other\textquotedblright\ means \textquotedblleft other than
$i_{1},i_{2},\ldots,i_{k}$\textquotedblright). In other words,%
\begin{equation}
\operatorname*{cyc}\nolimits_{i_{1},i_{2},\ldots,i_{k}}\left(  q\right)
=q\ \ \ \ \ \ \ \ \ \ \text{for every }q\in\left[  n\right]  \setminus\left\{
i_{1},i_{2},\ldots,i_{k}\right\}  . \label{sol.perm.c=ttt.cyc-manifest.2}%
\end{equation}

\end{itemize}

In the following, let $i_{k+1}$ mean $i_{1}$. Thus, $i_{k+1}=i_{1}$.

Define $\alpha\in S_{n}$ by $\alpha=\operatorname*{cyc}\nolimits_{i_{1}%
,i_{2},\ldots,i_{k}}$. Define $\beta\in S_{n}$ by $\beta=t_{i_{1},i_{2}}\circ
t_{i_{2},i_{3}}\circ\cdots\circ t_{i_{k-1},i_{k}}$.

\begin{vershort}
For each $j\in\left\{  1,2,\ldots,k-1\right\}  $, we define a transposition
$g_{j}$ in $S_{n}$ by $g_{j}=t_{i_{j},i_{j+1}}$. This transposition $g_{j}$
belongs to $S_{n}$, and thus is a bijective map from $\left[  n\right]  $ to
$\left[  n\right]  $.
\end{vershort}

\begin{verlong}
Recall that $S_{n}$ is the set of all permutations of the set $\left\{
1,2,\ldots,n\right\}  $. In other words, $S_{n}$ is the set of all
permutations of the set $\left[  n\right]  $ (since $\left\{  1,2,\ldots
,n\right\}  =\left[  n\right]  $). Now, $\alpha\in S_{n}$. In other words,
$\alpha$ is a permutation of the set $\left[  n\right]  $ (since $S_{n}$ is
the set of all permutations of the set $\left[  n\right]  $). In other words,
$\alpha$ is a bijective map from $\left[  n\right]  $ to $\left[  n\right]  $.
The same argument (applied to $\beta$ instead of $\alpha$) shows that $\beta$
is a bijective map from $\left[  n\right]  $ to $\left[  n\right]  $.

For each $j\in\left\{  1,2,\ldots,k-1\right\}  $, the transposition
$t_{i_{j},i_{j+1}}$ is a bijective map from $\left[  n\right]  $ to $\left[
n\right]  $\ \ \ \ \footnote{\textit{Proof.} Let $j\in\left\{  1,2,\ldots
,k-1\right\}  $. The transposition $t_{i_{j},i_{j+1}}$ belongs to $S_{n}$. In
other words, the transposition $t_{i_{j},i_{j+1}}$ is a permutation of the set
$\left[  n\right]  $ (since $S_{n}$ is the set of all permutations of the set
$\left[  n\right]  $). In other words, the transposition $t_{i_{j},i_{j+1}}$
is a bijective map from $\left[  n\right]  $ to $\left[  n\right]  $. Qed.}.
Thus, for each $j\in\left\{  1,2,\ldots,k-1\right\}  $, we can define a
bijective map $g_{j}$ from $\left[  n\right]  $ to $\left[  n\right]  $ by
$g_{j}=t_{i_{j},i_{j+1}}$. Consider these bijective maps $g_{j}$ for all
$j\in\left\{  1,2,\ldots,k-1\right\}  $.
\end{verlong}

For each $j\in\left\{  1,2,\ldots,k-1\right\}  $, we have $k-j\in\left\{
1,2,\ldots,k-1\right\}  $, and therefore $g_{k-j}$ is a well-defined bijective
map from $\left[  n\right]  $ to $\left[  n\right]  $. Hence, for each
$j\in\left\{  1,2,\ldots,k-1\right\}  $, we can define a bijective map $f_{j}$
from $\left[  n\right]  $ to $\left[  n\right]  $ by $f_{j}=g_{k-j}$. Consider
these bijective maps $f_{j}$ for all $j\in\left\{  1,2,\ldots,k-1\right\}  $.
We thus have defined $k-1$ bijective maps $f_{1},f_{2},\ldots,f_{k-1}$ from
$\left[  n\right]  $ to $\left[  n\right]  $. Note that%
\begin{align}
f_{k-1}\circ f_{\left(  k-1\right)  -1}\circ\cdots\circ f_{1}  &
=f_{k-1}\circ f_{k-2}\circ\cdots\circ f_{1}=\underbrace{g_{k-\left(
k-1\right)  }}_{=g_{1}}\circ\underbrace{g_{k-\left(  k-2\right)  }}_{=g_{2}%
}\circ\cdots\circ g_{k-1}\nonumber\\
&  \ \ \ \ \ \ \ \ \ \ \left(  \text{since }f_{j}=g_{k-j}\text{ for each }%
j\in\left\{  1,2,\ldots,k-1\right\}  \right) \nonumber\\
&  =g_{1}\circ g_{2}\circ\cdots\circ g_{k-1}=t_{i_{1},i_{2}}\circ
t_{i_{2},i_{3}}\circ\cdots\circ t_{i_{k-1},i_{k}}\nonumber\\
&  \ \ \ \ \ \ \ \ \ \ \left(  \text{since }g_{j}=t_{i_{j},i_{j+1}}\text{ for
each }j\in\left\{  1,2,\ldots,k-1\right\}  \right) \nonumber\\
&  =\beta\label{sol.perm.c=ttt.=beta}%
\end{align}
(since $\beta=t_{i_{1},i_{2}}\circ t_{i_{2},i_{3}}\circ\cdots\circ
t_{i_{k-1},i_{k}}$).

It is easy to see that each $j\in\left\{  1,2,\ldots,k-1\right\}  $ satisfies%
\begin{align}
&  f_{j}\left(  i_{k-j}\right)  =i_{k-j+1}\ \ \ \ \ \ \ \ \ \ \text{and}%
\label{sol.perm.c=ttt.fj1}\\
&  f_{j}\left(  i_{k-j+1}\right)  =i_{k-j}\ \ \ \ \ \ \ \ \ \ \text{and}%
\label{sol.perm.c=ttt.fj2}\\
&  \left(  f_{j}\left(  q\right)  =q\ \ \ \ \ \ \ \ \ \ \text{for each }%
q\in\left[  n\right]  \setminus\left\{  i_{k-j},i_{k-j+1}\right\}  \right)  .
\label{sol.perm.c=ttt.fj3}%
\end{align}

[\textit{Proof of (\ref{sol.perm.c=ttt.fj1}), (\ref{sol.perm.c=ttt.fj2}) and
(\ref{sol.perm.c=ttt.fj3}):} Let $j\in\left\{  1,2,\ldots,k-1\right\}  $. We
want to prove the three claims (\ref{sol.perm.c=ttt.fj1}),
(\ref{sol.perm.c=ttt.fj2}) and (\ref{sol.perm.c=ttt.fj3}).

We have $f_{j}=g_{k-j}$ (by the definition of $f_{j}$). Moreover,
$g_{k-j}=t_{i_{k-j},i_{\left(  k-j\right)  +1}}$ (by the definition of
$g_{k-j}$). Thus, $f_{j}=g_{k-j}=t_{i_{k-j},i_{\left(  k-j\right)  +1}%
}=t_{i_{k-j},i_{k-j+1}}$ (since $\left(  k-j\right)  +1=k-j+1$).

\begin{verlong}
From $j\in\left\{  1,2,\ldots,k-1\right\}  $, we obtain $k-j\in\left\{
1,2,\ldots,k-1\right\}  \subseteq\left\{  1,2,\ldots,k\right\}  $. Also,
\begin{align*}
k-j+1  &  =\left(  k-j\right)  +1\in\left\{  2,3,\ldots,k\right\}
\ \ \ \ \ \ \ \ \ \ \left(  \text{since }k-j\in\left\{  1,2,\ldots
,k-1\right\}  \right) \\
&  \subseteq\left\{  1,2,\ldots,k\right\}  .
\end{align*}
Thus, $k-j$ and $k-j+1$ are two elements of $\left\{  1,2,\ldots,k\right\}  $.
Thus, $i_{k-j}$ and $i_{k-j+1}$ are two elements of $\left[  n\right]  $
(since $i_{1},i_{2},\ldots,i_{k}$ are $k$ elements of $\left[  n\right]  $).
Moreover, $k-j\neq\left(  k-j\right)  +1=k-j+1$ and thus $i_{k-j}\neq
i_{k-j+1}$ (since the $k$ elements $i_{1},i_{2},\ldots,i_{k}$ are distinct).
In other words, $i_{k-j}$ and $i_{k-j+1}$ are distinct. Hence, $i_{k-j}$ and
$i_{k-j+1}$ are two distinct elements of $\left[  n\right]  $. In other words,
$i_{k-j}$ and $i_{k-j+1}$ are two distinct elements of $\left\{
1,2,\ldots,n\right\}  $ (since $\left[  n\right]  =\left\{  1,2,\ldots
,n\right\}  $). Thus, the transposition $t_{i_{k-j},i_{k-j+1}}$ in $S_{n}$ is
well-defined (according to Definition \ref{def.transpos}).
\end{verlong}

The definition of this transposition $t_{i_{k-j},i_{k-j+1}}$ shows that
$t_{i_{k-j},i_{k-j+1}}$ is the permutation in $S_{n}$ which swaps $i_{k-j}$
with $i_{k-j+1}$ while leaving all other elements of $\left\{  1,2,\ldots
,n\right\}  $ unchanged. In other words, $f_{j}$ is the permutation in $S_{n}$
which swaps $i_{k-j}$ with $i_{k-j+1}$ while leaving all other elements of
$\left[  n\right]  $ unchanged (since $f_{j}=t_{i_{k-j},i_{k-j+1}}$ and
$\left[  n\right]  =\left\{  1,2,\ldots,n\right\}  $). In other words, $f_{j}$
is the permutation in $S_{n}$ that satisfies $f_{j}\left(  i_{k-j}\right)
=i_{k-j+1}$ and $f_{j}\left(  i_{k-j+1}\right)  =i_{k-j}$ and%
\[
\left(  f_{j}\left(  q\right)  =q\ \ \ \ \ \ \ \ \ \ \text{for each }%
q\in\left[  n\right]  \setminus\left\{  i_{k-j},i_{k-j+1}\right\}  \right)  .
\]
Thus, the three claims (\ref{sol.perm.c=ttt.fj1}), (\ref{sol.perm.c=ttt.fj2})
and (\ref{sol.perm.c=ttt.fj3}) are proven.]

Now, we claim that each $x\in\left[  n\right]  $ satisfies%
\begin{equation}
\alpha\left(  x\right)  =\beta\left(  x\right)  . \label{sol.perm.c=ttt.aq=bq}%
\end{equation}

[\textit{Proof of (\ref{sol.perm.c=ttt.aq=bq}):} Let $x\in\left[  n\right]  $.
We are in one of the following three cases:

\textit{Case 1:} We have $x\notin\left\{  i_{1},i_{2},\ldots,i_{k}\right\}  $.

\textit{Case 2:} We have $x=i_{k}$.

\textit{Case 3:} We have neither $x\notin\left\{  i_{1},i_{2},\ldots
,i_{k}\right\}  $ nor $x=i_{k}$.

\begin{vershort}
Let us first consider Case 1. In this case, we have $x\notin\left\{
i_{1},i_{2},\ldots,i_{k}\right\}  $. Combining $x\in\left[  n\right]  $ with
$x\notin\left\{  i_{1},i_{2},\ldots,i_{k}\right\}  $, we obtain $x\in\left[
n\right]  \setminus\left\{  i_{1},i_{2},\ldots,i_{k}\right\}  $. Hence,
(\ref{sol.perm.c=ttt.cyc-manifest.2}) (applied to $q=x$) yields
$\operatorname*{cyc}\nolimits_{i_{1},i_{2},\ldots,i_{k}}\left(  x\right)  =x$.
On the other hand, we have $f_{j}\left(  x\right)  =x$ for each $j\in\left\{
1,2,\ldots,k-1\right\}  $\ \ \ \ \footnote{\textit{Proof.} Let $j\in\left\{
1,2,\ldots,k-1\right\}  $. Thus, $x\in\left[  n\right]  \setminus
\underbrace{\left\{  i_{1},i_{2},\ldots,i_{k}\right\}  }_{\supseteq\left\{
i_{k-j},i_{k-j+1}\right\}  }\subseteq\left[  n\right]  \setminus\left\{
i_{k-j},i_{k-j+1}\right\}  $. Hence, (\ref{sol.perm.c=ttt.fj3}) (applied to
$q=x$) yields $f_{j}\left(  x\right)  =x$. Qed.}. Hence, Lemma
\ref{lem.sol.sigmacrosstau.compose0} (applied to $m=k-1$ and $X=\left[
n\right]  $) yields $\left(  f_{k-1}\circ f_{\left(  k-1\right)  -1}%
\circ\cdots\circ f_{1}\right)  \left(  x\right)  =x$. In view of
(\ref{sol.perm.c=ttt.=beta}), this rewrites as $\beta\left(  x\right)  =x$.
Comparing this with $\underbrace{\alpha}_{=\operatorname*{cyc}\nolimits_{i_{1}%
,i_{2},\ldots,i_{k}}}\left(  x\right)  =\operatorname*{cyc}\nolimits_{i_{1}%
,i_{2},\ldots,i_{k}}\left(  x\right)  =x$, we obtain $\alpha\left(  x\right)
=\beta\left(  x\right)  $. Hence, (\ref{sol.perm.c=ttt.aq=bq}) is proven in
Case 1.
\end{vershort}

\begin{verlong}
Let us first consider Case 1. In this case, we have $x\notin\left\{
i_{1},i_{2},\ldots,i_{k}\right\}  $. Combining $x\in\left[  n\right]  $ with
$x\notin\left\{  i_{1},i_{2},\ldots,i_{k}\right\}  $, we obtain $x\in\left[
n\right]  \setminus\left\{  i_{1},i_{2},\ldots,i_{k}\right\}  $. Hence,
(\ref{sol.perm.c=ttt.cyc-manifest.2}) (applied to $q=x$) yields
$\operatorname*{cyc}\nolimits_{i_{1},i_{2},\ldots,i_{k}}\left(  x\right)  =x$.
On the other hand, we have $f_{j}\left(  x\right)  =x$ for each $j\in\left\{
1,2,\ldots,k-1\right\}  $\ \ \ \ \footnote{\textit{Proof.} Let $j\in\left\{
1,2,\ldots,k-1\right\}  $. Thus, $k-j\in\left\{  1,2,\ldots,k-1\right\}
\subseteq\left\{  1,2,\ldots,k\right\}  $ and thus $i_{k-j}\in\left\{
i_{1},i_{2},\ldots,i_{k}\right\}  $. If we had $x=i_{k-j}$, then we would have
$x=i_{k-j}\in\left\{  i_{1},i_{2},\ldots,i_{k}\right\}  $, which would
contradict $x\notin\left\{  i_{1},i_{2},\ldots,i_{k}\right\}  $. Hence, we
cannot have $x=i_{k-j}$. In other words, we have $x\neq i_{k-j}$.
\par
Also,
\begin{align*}
k-j+1  &  =\left(  k-j\right)  +1\in\left\{  2,3,\ldots,k\right\}
\ \ \ \ \ \ \ \ \ \ \left(  \text{since }k-j\in\left\{  1,2,\ldots
,k-1\right\}  \right) \\
&  \subseteq\left\{  1,2,\ldots,k\right\}
\end{align*}
and thus $i_{k-j+1}\in\left\{  i_{1},i_{2},\ldots,i_{k}\right\}  $. If we had
$x=i_{k-j+1}$, then we would have $x=i_{k-j+1}\in\left\{  i_{1},i_{2}%
,\ldots,i_{k}\right\}  $, which would contradict $x\notin\left\{  i_{1}%
,i_{2},\ldots,i_{k}\right\}  $. Hence, we cannot have $x=i_{k-j+1}$. In other
words, we have $x\neq i_{k-j+1}$.
\par
We have $\left(  \text{neither }x=i_{k-j}\text{ nor }x=i_{k-j+1}\right)  $
(since $x\neq i_{k-j}$ and $x\neq i_{k-j+1}$).
\par
If we had $x\in\left\{  i_{k-j},i_{k-j+1}\right\}  $, then we would have
$\left(  \text{either }x=i_{k-j}\text{ or }x=i_{k-j+1}\right)  $, which would
contradict $\left(  \text{neither }x=i_{k-j}\text{ nor }x=i_{k-j+1}\right)  $.
Thus, we cannot have $x\in\left\{  i_{k-j},i_{k-j+1}\right\}  $. Hence, we
have $x\notin\left\{  i_{k-j},i_{k-j+1}\right\}  $. Combining $x\in\left[
n\right]  $ with $x\notin\left\{  i_{k-j},i_{k-j+1}\right\}  $, we obtain
$x\in\left[  n\right]  \setminus\left\{  i_{k-j},i_{k-j+1}\right\}  $. Hence,
(\ref{sol.perm.c=ttt.fj3}) (applied to $q=x$) yields $f_{j}\left(  x\right)
=x$. Qed.}. Hence, Lemma \ref{lem.sol.sigmacrosstau.compose0} (applied to
$m=k-1$ and $X=\left[  n\right]  $) yields $\left(  f_{k-1}\circ f_{\left(
k-1\right)  -1}\circ\cdots\circ f_{1}\right)  \left(  x\right)  =x$. In view
of (\ref{sol.perm.c=ttt.=beta}), this rewrites as $\beta\left(  x\right)  =x$.
Comparing this with%
\[
\underbrace{\alpha}_{=\operatorname*{cyc}\nolimits_{i_{1},i_{2},\ldots,i_{k}}%
}\left(  x\right)  =\operatorname*{cyc}\nolimits_{i_{1},i_{2},\ldots,i_{k}%
}\left(  x\right)  =x,
\]
we obtain $\alpha\left(  x\right)  =\beta\left(  x\right)  $. Hence,
(\ref{sol.perm.c=ttt.aq=bq}) is proven in Case 1.
\end{verlong}

Let us next consider Case 2. In this case, we have $x=i_{k}$.

\begin{verlong}
Recall that $i_{k},i_{k-1},\ldots,i_{1}$ are $k$ elements of $\left[
n\right]  $. In other words, $i_{k-0},i_{k-1},\ldots,i_{k-\left(  k-1\right)
}$ are $k$ elements of $\left[  n\right]  $ (since the elements $i_{k-0}%
,i_{k-1},\ldots,i_{k-\left(  k-1\right)  }$ are precisely the elements
$i_{k},i_{k-1},\ldots,i_{1}$). In other words, $i_{k-0},i_{k-1},\ldots
,i_{k-\left(  k-1\right)  }$ are $\left(  k-1\right)  +1$ elements of $\left[
n\right]  $ (since $\left(  k-1\right)  +1=k$).
\end{verlong}

We have $f_{j}\left(  i_{k-\left(  j-1\right)  }\right)  =i_{k-j}$ for each
$j\in\left\{  1,2,\ldots,k-1\right\}  $\ \ \ \ \footnote{\textit{Proof.} Let
$j\in\left\{  1,2,\ldots,k-1\right\}  $. Then, (\ref{sol.perm.c=ttt.fj2})
yields $f_{j}\left(  i_{k-j+1}\right)  =i_{k-j}$. In view of $k-j+1=k-\left(
j-1\right)  $, this rewrites as $f_{j}\left(  i_{k-\left(  j-1\right)
}\right)  =i_{k-j}$. Qed.}. Thus, Lemma \ref{lem.sol.sigmacrosstau.composeall}
(applied to $m=k-1$ and $X=\left[  n\right]  $ and $x_{j}=i_{k-j}$) yields
\newline$\left(  f_{k-1}\circ f_{\left(  k-1\right)  -1}\circ\cdots\circ
f_{1}\right)  \left(  i_{k-0}\right)  =i_{k-\left(  k-1\right)  }$. In view of
$k-0=k$ and $k-\left(  k-1\right)  =1$, this rewrites as $\left(  f_{k-1}\circ
f_{\left(  k-1\right)  -1}\circ\cdots\circ f_{1}\right)  \left(  i_{k}\right)
=i_{1}$. In view of (\ref{sol.perm.c=ttt.=beta}), this rewrites as
$\beta\left(  i_{k}\right)  =i_{1}$.

On the other hand, $k\in\left\{  1,2,\ldots,k\right\}  $ (since $k\geq1$).
Hence, (\ref{sol.perm.c=ttt.cyc-manifest.1}) (applied to $p=k$) yields
$\operatorname*{cyc}\nolimits_{i_{1},i_{2},\ldots,i_{k}}\left(  i_{k}\right)
=i_{k+1}=i_{1}$. Thus,%
\[
\underbrace{\alpha}_{=\operatorname*{cyc}\nolimits_{i_{1},i_{2},\ldots,i_{k}}%
}\left(  \underbrace{x}_{=i_{k}}\right)  =\operatorname*{cyc}\nolimits_{i_{1}%
,i_{2},\ldots,i_{k}}\left(  i_{k}\right)  =i_{1}.
\]
Comparing this with $\beta\left(  \underbrace{x}_{=i_{k}}\right)
=\beta\left(  i_{k}\right)  =i_{1}$, we obtain $\alpha\left(  x\right)
=\beta\left(  x\right)  $. Hence, (\ref{sol.perm.c=ttt.aq=bq}) is proven in
Case 2.

Let us finally consider Case 3. In this case, we have neither $x\notin\left\{
i_{1},i_{2},\ldots,i_{k}\right\}  $ nor $x=i_{k}$. Thus, we don't have
$x\notin\left\{  i_{1},i_{2},\ldots,i_{k}\right\}  $. Hence, we have
$x\in\left\{  i_{1},i_{2},\ldots,i_{k}\right\}  $. In other words, $x=i_{p}$
for some $p\in\left\{  1,2,\ldots,k\right\}  $. Consider this $p$.

Also, we don't have $x=i_{k}$ (since we have neither $x\notin\left\{
i_{1},i_{2},\ldots,i_{k}\right\}  $ nor $x=i_{k}$). In other words, we have
$x\neq i_{k}$. From $x=i_{p}$, we obtain $i_{p}=x\neq i_{k}$, thus $p\neq k$.
Combining $p\in\left\{  1,2,\ldots,k\right\}  $ with $p\neq k$, we obtain
$p\in\left\{  1,2,\ldots,k\right\}  \setminus\left\{  k\right\}  =\left\{
1,2,\ldots,k-1\right\}  $. Hence, $k-p\in\left\{  1,2,\ldots,k-1\right\}  $.
Thus, the map $f_{k-p}$ is well-defined.

We have $p\in\left\{  1,2,\ldots,k-1\right\}  $, thus $p+1\in\left\{
2,3,\ldots,k\right\}  \subseteq\left\{  1,2,\ldots,k\right\}  $. Hence,
$i_{p+1}\in\left\{  i_{1},i_{2},\ldots,i_{k}\right\}  \subseteq\left[
n\right]  $ (since $i_{1},i_{2},\ldots,i_{k}$ are elements of $\left[
n\right]  $). Thus, we can define $y\in\left[  n\right]  $ by $y=i_{p+1}$.
Consider this $y$. Thus,%
\[
x=i_{p}\ \ \ \ \ \ \ \ \ \ \text{and}\ \ \ \ \ \ \ \ \ \ y=i_{p+1}.
\]

The equality (\ref{sol.perm.c=ttt.fj1}) (applied to $j=k-p$) yields
$f_{k-p}\left(  i_{k-\left(  k-p\right)  }\right)  =i_{k-\left(  k-p\right)
+1}$. In view of $k-\left(  k-p\right)  =p$, this rewrites as $f_{k-p}\left(
i_{p}\right)  =i_{p+1}$. In view of $i_{p}=x$ and $i_{p+1}=y$, this rewrites
as
\begin{equation}
f_{k-p}\left(  x\right)  =y. \label{sol.perm.c=ttt.aq=bq.c3.1}%
\end{equation}

Moreover, we have%
\begin{equation}
f_{j}\left(  x\right)  =x\ \ \ \ \ \ \ \ \ \ \text{for each }j\in\left\{
1,2,\ldots,k-1\right\}  \text{ satisfying }j<k-p.
\label{sol.perm.c=ttt.aq=bq.c3.2}%
\end{equation}

\begin{vershort}
[\textit{Proof of (\ref{sol.perm.c=ttt.aq=bq.c3.2}):} Let $j\in\left\{
1,2,\ldots,k-1\right\}  $ be such that $j<k-p$.

From $j<k-p$, we obtain $k-j>p$, thus $k-j\neq p$. Hence, $i_{k-j}\neq i_{p}$
(since the $k$ elements $i_{1},i_{2},\ldots,i_{k}$ are distinct). In other
words, $i_{p}\neq i_{k-j}$. Hence, $x=i_{p}\neq i_{k-j}$.

We have $k-j+1>k-j>p$, thus $k-j+1\neq p$. Hence, $i_{k-j+1}\neq i_{p}$ (since
the $k$ elements $i_{1},i_{2},\ldots,i_{k}$ are distinct). In other words,
$i_{p}\neq i_{k-j+1}$. Hence, $x=i_{p}\neq i_{k-j+1}$.

Combining $x\neq i_{k-j}$ with $x\neq i_{k-j+1}$, we obtain $x\notin\left\{
i_{k-j},i_{k-j+1}\right\}  $. Combining $x\in\left[  n\right]  $ with this, we
obtain $x\in\left[  n\right]  \setminus\left\{  i_{k-j},i_{k-j+1}\right\}  $.
Hence, (\ref{sol.perm.c=ttt.fj3}) (applied to $q=x$) yields $f_{j}\left(
x\right)  =x$. This proves (\ref{sol.perm.c=ttt.aq=bq.c3.2}).]
\end{vershort}

\begin{verlong}
[\textit{Proof of (\ref{sol.perm.c=ttt.aq=bq.c3.2}):} Let $j\in\left\{
1,2,\ldots,k-1\right\}  $ be such that $j<k-p$. We must prove that
$f_{j}\left(  x\right)  =x$.

From $j\in\left\{  1,2,\ldots,k-1\right\}  $, we obtain $k-j\in\left\{
1,2,\ldots,k-1\right\}  \subseteq\left\{  1,2,\ldots,k\right\}  $. Hence,
$i_{k-j}$ is well-defined. Also,%
\begin{align*}
k-j+1  &  =\left(  k-j\right)  +1\in\left\{  2,3,\ldots,k\right\}
\ \ \ \ \ \ \ \ \ \ \left(  \text{since }k-j\in\left\{  1,2,\ldots
,k-1\right\}  \right) \\
&  \subseteq\left\{  1,2,\ldots,k\right\}  .
\end{align*}
Hence, $i_{k-j+1}$ is well-defined.

We have $k-\underbrace{j}_{<k-p}>k-\left(  k-p\right)  =p$, thus $k-j\neq p$.
Hence, $i_{k-j}\neq i_{p}$ (since the $k$ elements $i_{1},i_{2},\ldots,i_{k}$
are distinct). In other words, $i_{p}\neq i_{k-j}$. Hence, $x=i_{p}\neq
i_{k-j}$.

We have $k-j+1>k-j>p$, thus $k-j+1\neq p$. Hence, $i_{k-j+1}\neq i_{p}$ (since
the $k$ elements $i_{1},i_{2},\ldots,i_{k}$ are distinct). In other words,
$i_{p}\neq i_{k-j+1}$. Hence, $x=i_{p}\neq i_{k-j+1}$.

We have $\left(  \text{neither }x=i_{k-j}\text{ nor }x=i_{k-j+1}\right)  $
(since $x\neq i_{k-j}$ and $x\neq i_{k-j+1}$). If we had $x\in\left\{
i_{k-j},i_{k-j+1}\right\}  $, then we would have $\left(  \text{either
}x=i_{k-j}\text{ or }x=i_{k-j+1}\right)  $, which would contradict $\left(
\text{neither }x=i_{k-j}\text{ nor }x=i_{k-j+1}\right)  $. Thus, we cannot
have $x\in\left\{  i_{k-j},i_{k-j+1}\right\}  $. Hence, we have $x\notin%
\left\{  i_{k-j},i_{k-j+1}\right\}  $. Combining $x\in\left[  n\right]  $ with
$x\notin\left\{  i_{k-j},i_{k-j+1}\right\}  $, we obtain $x\in\left[
n\right]  \setminus\left\{  i_{k-j},i_{k-j+1}\right\}  $. Hence,
(\ref{sol.perm.c=ttt.fj3}) (applied to $q=x$) yields $f_{j}\left(  x\right)
=x$. This proves (\ref{sol.perm.c=ttt.aq=bq.c3.2}).]
\end{verlong}

Furthermore, we have%
\begin{equation}
f_{j}\left(  y\right)  =y\ \ \ \ \ \ \ \ \ \ \text{for each }j\in\left\{
1,2,\ldots,k-1\right\}  \text{ satisfying }j>k-p.
\label{sol.perm.c=ttt.aq=bq.c3.3}%
\end{equation}

\begin{vershort}
[\textit{Proof of (\ref{sol.perm.c=ttt.aq=bq.c3.3}):} Let $j\in\left\{
1,2,\ldots,k-1\right\}  $ be such that $j>k-p$.

From $j>k-p$, we obtain $k-j<p<p+1$, thus $k-j\neq p+1$. Hence, $i_{k-j}\neq
i_{p+1}$ (since the $k$ elements $i_{1},i_{2},\ldots,i_{k}$ are distinct). In
other words, $i_{p+1}\neq i_{k-j}$. Hence, $y=i_{p+1}\neq i_{k-j}$.

We have $k-j<p$, thus $k-j\neq p$ and therefore $k-j+1\neq p+1$. Hence,
$i_{k-j+1}\neq i_{p+1}$ (since the $k$ elements $i_{1},i_{2},\ldots,i_{k}$ are
distinct). In other words, $i_{p+1}\neq i_{k-j+1}$. Hence, $y=i_{p+1}\neq
i_{k-j+1}$.

Combining $y\neq i_{k-j}$ with $y\neq i_{k-j+1}$, we obtain $y\notin\left\{
i_{k-j},i_{k-j+1}\right\}  $. Combining $y\in\left[  n\right]  $ with this, we
obtain $y\in\left[  n\right]  \setminus\left\{  i_{k-j},i_{k-j+1}\right\}  $.
Hence, (\ref{sol.perm.c=ttt.fj3}) (applied to $q=y$) yields $f_{j}\left(
y\right)  =y$. This proves (\ref{sol.perm.c=ttt.aq=bq.c3.3}).]
\end{vershort}

\begin{verlong}
[\textit{Proof of (\ref{sol.perm.c=ttt.aq=bq.c3.3}):} Let $j\in\left\{
1,2,\ldots,k-1\right\}  $ be such that $j>k-p$. We must prove that
$f_{j}\left(  y\right)  =y$.

From $j\in\left\{  1,2,\ldots,k-1\right\}  $, we obtain $k-j\in\left\{
1,2,\ldots,k-1\right\}  \subseteq\left\{  1,2,\ldots,k\right\}  $. Hence,
$i_{k-j}$ is well-defined. Also,%
\begin{align*}
k-j+1  &  =\left(  k-j\right)  +1\in\left\{  2,3,\ldots,k\right\}
\ \ \ \ \ \ \ \ \ \ \left(  \text{since }k-j\in\left\{  1,2,\ldots
,k-1\right\}  \right) \\
&  \subseteq\left\{  1,2,\ldots,k\right\}  .
\end{align*}
Hence, $i_{k-j+1}$ is well-defined.

We have $k-\underbrace{j}_{>k-p}<k-\left(  k-p\right)  =p$, thus
$\underbrace{k-j}_{<p}+1<p+1$ and therefore $k-j+1\neq p+1$. Hence,
$i_{k-j+1}\neq i_{p+1}$ (since the $k$ elements $i_{1},i_{2},\ldots,i_{k}$ are
distinct). In other words, $i_{p+1}\neq i_{k-j+1}$. Hence, $y=i_{p+1}\neq
i_{k-j+1}$.

We have $k-j<p<p+1$, thus $k-j\neq p+1$. Hence, $i_{k-j}\neq i_{p+1}$ (since
the $k$ elements $i_{1},i_{2},\ldots,i_{k}$ are distinct). In other words,
$i_{p+1}\neq i_{k-j}$. Hence, $y=i_{p+1}\neq i_{k-j}$.

We have $\left(  \text{neither }y=i_{k-j}\text{ nor }y=i_{k-j+1}\right)  $
(since $y\neq i_{k-j}$ and $y\neq i_{k-j+1}$). If we had $y\in\left\{
i_{k-j},i_{k-j+1}\right\}  $, then we would have $\left(  \text{either
}y=i_{k-j}\text{ or }y=i_{k-j+1}\right)  $, which would contradict $\left(
\text{neither }y=i_{k-j}\text{ nor }y=i_{k-j+1}\right)  $. Thus, we cannot
have $y\in\left\{  i_{k-j},i_{k-j+1}\right\}  $. Hence, we have $y\notin%
\left\{  i_{k-j},i_{k-j+1}\right\}  $. Combining $y\in\left[  n\right]  $ with
$y\notin\left\{  i_{k-j},i_{k-j+1}\right\}  $, we obtain $y\in\left[
n\right]  \setminus\left\{  i_{k-j},i_{k-j+1}\right\}  $. Hence,
(\ref{sol.perm.c=ttt.fj3}) (applied to $q=y$) yields $f_{j}\left(  y\right)
=y$. This proves (\ref{sol.perm.c=ttt.aq=bq.c3.3}).]
\end{verlong}

Now, we have proven the equalities (\ref{sol.perm.c=ttt.aq=bq.c3.1}),
(\ref{sol.perm.c=ttt.aq=bq.c3.2}) and (\ref{sol.perm.c=ttt.aq=bq.c3.3}).
Hence, Lemma \ref{lem.sol.sigmacrosstau.compose1} (applied to $k-1$, $\left[
n\right]  $ and $k-p$ instead of $m$, $X$ and $i$) yields%
\[
\left(  f_{k-1}\circ f_{\left(  k-1\right)  -1}\circ\cdots\circ f_{1}\right)
\left(  x\right)  =y.
\]
In view of (\ref{sol.perm.c=ttt.=beta}), this rewrites as $\beta\left(
x\right)  =y$. Comparing this with%
\begin{align*}
\underbrace{\alpha}_{=\operatorname*{cyc}\nolimits_{i_{1},i_{2},\ldots,i_{k}}%
}\left(  \underbrace{x}_{=i_{p}}\right)   &  =\operatorname*{cyc}%
\nolimits_{i_{1},i_{2},\ldots,i_{k}}\left(  i_{p}\right)  =i_{p+1}%
\ \ \ \ \ \ \ \ \ \ \left(  \text{by (\ref{sol.perm.c=ttt.cyc-manifest.1}%
)}\right) \\
&  =y,
\end{align*}
we obtain $\alpha\left(  x\right)  =\beta\left(  x\right)  $. Hence,
(\ref{sol.perm.c=ttt.aq=bq}) is proven in Case 3.

We have now proven (\ref{sol.perm.c=ttt.aq=bq}) in each of the three Cases 1,
2 and 3. Since these three Cases cover all possibilities, we thus conclude
that (\ref{sol.perm.c=ttt.aq=bq}) always holds.]

But from (\ref{sol.perm.c=ttt.aq=bq}), we immediately obtain $\alpha=\beta$
(since both $\alpha$ and $\beta$ are maps from $\left[  n\right]  $ to
$\left[  n\right]  $). In view of $\alpha=\operatorname*{cyc}\nolimits_{i_{1}%
,i_{2},\ldots,i_{k}}$ and $\beta=t_{i_{1},i_{2}}\circ t_{i_{2},i_{3}}%
\circ\cdots\circ t_{i_{k-1},i_{k}}$, this rewrites as
\[
\operatorname*{cyc}\nolimits_{i_{1},i_{2},\ldots,i_{k}}=t_{i_{1},i_{2}}\circ
t_{i_{2},i_{3}}\circ\cdots\circ t_{i_{k-1},i_{k}}.
\]
This solves Exercise \ref{exe.perm.c=ttt}.
\end{proof}

\subsubsection{A particular case}

In this subsection, we shall derive a particular case of Exercise
\ref{exe.perm.c=ttt} that will come useful later. We make a definition:

\begin{definition}
\label{def.sol.perm.c=ttt.cuv}Let $n\in\mathbb{N}$. Let $u$ and $v$ be two
elements of $\left[  n\right]  $ such that $u\leq v$. Then, we define a
permutation $c_{u,v}\in S_{n}$ by%
\[
c_{u,v}=\operatorname*{cyc}\nolimits_{v,v-1,v-2,\ldots,u}.
\]
(This is well-defined by Corollary \ref{cor.sol.perm.c=ttt.cuv1} \textbf{(a)} below.)
\end{definition}

\begin{corollary}
\label{cor.sol.perm.c=ttt.cuv1}Let $n\in\mathbb{N}$. Let $u$ and $v$ be two
elements of $\left[  n\right]  $ such that $u\leq v$. Then:

\textbf{(a)} The permutation $c_{u,v}$ is well-defined.

\textbf{(b)} We have $c_{u,v}=s_{v-1}\circ s_{v-2}\circ\cdots\circ s_{u}$.
\end{corollary}

Before we prove this corollary, let us establish an almost trivial lemma:

\begin{lemma}
\label{lem.sol.perm.c=ttt.s=t}Let $n\in\mathbb{N}$. Let $i\in\left\{
1,2,\ldots,n-1\right\}  $. Then, $t_{i+1,i}=s_{i}$.
\end{lemma}

\begin{proof}
[Proof of Lemma \ref{lem.sol.perm.c=ttt.s=t}.]Recall that $s_{i}$ is the
permutation in $S_{n}$ that swaps $i$ with $i+1$ but leaves all other numbers
unchanged (by the definition of $s_{i}$).

\begin{vershort}
On the other hand, $i\in\left\{  1,2,\ldots,n-1\right\}  $. Hence, $i+1$ and
$i$ are two distinct elements of $\left\{  1,2,\ldots,n\right\}  $. Thus,
$t_{i+1,i}$ is the permutation in $S_{n}$ which swaps $i+1$ with $i$ while
leaving all other elements of $\left\{  1,2,\ldots,n\right\}  $ unchanged (by
the definition of $t_{i+1,i}$). In other words, $t_{i+1,i}$ is the permutation
in $S_{n}$ which swaps $i$ with $i+1$ while leaving all other elements of
$\left\{  1,2,\ldots,n\right\}  $ unchanged (because swapping $i+1$ with $i$
is tantamount to swapping $i$ with $i+1$). In other words, $t_{i+1,i}$ is the
permutation in $S_{n}$ that swaps $i$ with $i+1$ but leaves all other numbers unchanged.
\end{vershort}

\begin{verlong}
On the other hand, $i\in\left\{  1,2,\ldots,n-1\right\}  \subseteq\left\{
1,2,\ldots,n\right\}  $ and
\begin{align*}
i+1  &  \in\left\{  2,3,\ldots,n\right\}  \ \ \ \ \ \ \ \ \ \ \left(
\text{since }i\in\left\{  1,2,\ldots,n-1\right\}  \right) \\
&  \subseteq\left\{  1,2,\ldots,n\right\}  .
\end{align*}
Hence, $i+1$ and $i$ are two elements of $\left\{  1,2,\ldots,n\right\}  $.
These two elements $i+1$ and $i$ are distinct (since $i+1\neq i$). Hence,
$t_{i+1,i}$ is the permutation in $S_{n}$ which swaps $i+1$ with $i$ while
leaving all other elements of $\left\{  1,2,\ldots,n\right\}  $ unchanged (by
the definition of $t_{i+1,i}$). In other words, $t_{i+1,i}$ is the permutation
in $S_{n}$ which swaps $i$ with $i+1$ while leaving all other elements of
$\left\{  1,2,\ldots,n\right\}  $ unchanged (because swapping $i+1$ with $i$
is tantamount to swapping $i$ with $i+1$). In other words, $t_{i+1,i}$ is the
permutation in $S_{n}$ that swaps $i$ with $i+1$ but leaves all other numbers unchanged.
\end{verlong}

Now, we have shown that both $t_{i+1,i}$ and $s_{i}$ are the permutation in
$S_{n}$ that swaps $i$ with $i+1$ but leaves all other numbers unchanged.
Hence, $t_{i+1,i}$ and $s_{i}$ are the same permutation. In other words,
$t_{i+1,i}=s_{i}$. This proves Lemma \ref{lem.sol.perm.c=ttt.s=t}.
\end{proof}

\begin{proof}
[Proof of Corollary \ref{cor.sol.perm.c=ttt.cuv1}.]We have $u\in\left[
n\right]  =\left\{  1,2,\ldots,n\right\}  $ (by the definition of $\left[
n\right]  $), thus $u\geq1$. Hence, $1\leq u$. Also, $v\in\left[  n\right]
=\left\{  1,2,\ldots,n\right\}  $, thus $v\leq n$. Hence, $1\leq u\leq v\leq
n$.

\begin{vershort}
Also, $\underbrace{v}_{\leq n}-\underbrace{u}_{\geq1}\leq n-1$. Combining this
with $v-\underbrace{u}_{\leq v}\geq v-v=0$, we obtain $v-u\in\left\{
0,1,\ldots,n-1\right\}  $, so that $v-u+1\in\left\{  1,2,\ldots,n\right\}  $.
\end{vershort}

\begin{verlong}
Also, $\underbrace{v}_{\leq n}-\underbrace{u}_{\geq1}\leq n-1$. Combining this
with $v-\underbrace{u}_{\leq v}\geq v-v=0$, we obtain $v-u\in\left\{
0,1,\ldots,n-1\right\}  $, so that $\left(  v-u\right)  +1\in\left\{
1,2,\ldots,n\right\}  $. Hence, $v-u+1=\left(  v-u\right)  +1\in\left\{
1,2,\ldots,n\right\}  $.
\end{verlong}

\begin{vershort}
Now, $v,v-1,v-2,\ldots,u$ are $v-u+1$ distinct elements of the set $\left\{
1,2,\ldots,n\right\}  $ (since $1\leq u\leq v\leq n$). Hence, the permutation
$\operatorname*{cyc}\nolimits_{v,v-1,v-2,\ldots,u}$ is well-defined. In other
words, the permutation $c_{u,v}$ is well-defined (since $c_{u,v}$ was defined
by $c_{u,v}=\operatorname*{cyc}\nolimits_{v,v-1,v-2,\ldots,u}$). This proves
Corollary \ref{cor.sol.perm.c=ttt.cuv1} \textbf{(a)}.
\end{vershort}

\begin{verlong}
From $u\leq v$, we conclude that $v,v-1,v-2,\ldots,u$ are $v-u+1$ numbers.
Also, both numbers $u$ and $v$ belong to $\left\{  1,2,\ldots,n\right\}  $
(since $u\in\left\{  1,2,\ldots,n\right\}  $ and $v\in\left\{  1,2,\ldots
,n\right\}  $); thus, all the numbers between $u$ and $v$ also belong to
$\left\{  1,2,\ldots,n\right\}  $. In other words, the $v-u+1$ numbers
$v,v-1,v-2,\ldots,u$ all belong to $\left\{  1,2,\ldots,n\right\}  $. Hence,
$v,v-1,v-2,\ldots,u$ are $v-u+1$ distinct elements of the set $\left\{
1,2,\ldots,n\right\}  $. Hence, the permutation $\operatorname*{cyc}%
\nolimits_{v,v-1,v-2,\ldots,u}$ is well-defined (according to Definition
\ref{def.perm.cycles}). In other words, the permutation $c_{u,v}$ is
well-defined (since $c_{u,v}$ was defined by $c_{u,v}=\operatorname*{cyc}%
\nolimits_{v,v-1,v-2,\ldots,u}$). This proves Corollary
\ref{cor.sol.perm.c=ttt.cuv1} \textbf{(a)}.
\end{verlong}

\begin{vershort}
\textbf{(b)} Each $i\in\left\{  v-1,v-2,\ldots,u\right\}  $ satisfies
$t_{i+1,i}=s_{i}$ (by Lemma \ref{lem.sol.perm.c=ttt.s=t}) and thus
$s_{i}=t_{i+1,i}$. In other words,%
\begin{align*}
\left(  s_{v-1},s_{v-2},\ldots,s_{u}\right)   &  =\left(  t_{\left(
v-1\right)  +1,v-1},t_{\left(  v-2\right)  +1,v-2},\ldots,t_{u+1,u}\right) \\
&  =\left(  t_{v,v-1},t_{v-1,v-2},\ldots,t_{u+1,u}\right)  .
\end{align*}
Hence,
\begin{equation}
s_{v-1}\circ s_{v-2}\circ\cdots\circ s_{u}=t_{v,v-1}\circ t_{v-1,v-2}%
\circ\cdots\circ t_{u+1,u}. \label{pf.cor.sol.perm.c=ttt.cuv1.b.short.1}%
\end{equation}

But Exercise \ref{exe.perm.c=ttt} (applied to $k=v-u+1$ and $\left(
i_{1},i_{2},\ldots,i_{k}\right)  =\left(  v,v-1,v-2,\ldots,u\right)  $) yields%
\[
\operatorname*{cyc}\nolimits_{v,v-1,v-2,\ldots,u}=t_{v,v-1}\circ
t_{v-1,v-2}\circ\cdots\circ t_{u+1,u}=s_{v-1}\circ s_{v-2}\circ\cdots\circ
s_{u}%
\]
(by (\ref{pf.cor.sol.perm.c=ttt.cuv1.b.short.1})). Now, the definition of
$c_{u,v}$ yields $c_{u,v}=\operatorname*{cyc}\nolimits_{v,v-1,v-2,\ldots
,u}=s_{v-1}\circ s_{v-2}\circ\cdots\circ s_{u}$. This proves Corollary
\ref{cor.sol.perm.c=ttt.cuv1} \textbf{(b)}. \qedhere

\end{vershort}

\begin{verlong}
\textbf{(b)} Let $k=v-u+1$. Thus, $k=v-u+1\in\left\{  1,2,\ldots,n\right\}  $
and $v-\underbrace{k}_{=v-u+1}+1=v-\left(  v-u+1\right)  +1=u$.

For each $p\in\left\{  1,2,\ldots,k\right\}  $, define the integer $i_{p}$ by
$i_{p}=v-p+1$. Thus,%
\begin{align*}
\left(  i_{1},i_{2},\ldots,i_{k}\right)   &  =\left(  v-1+1,v-2+1,\ldots
,v-k+1\right) \\
&  =\left(  v,v-1,v-2,\ldots,v-k+1\right)  =\left(  v,v-1,v-2,\ldots,u\right)
\end{align*}
(since $v-k+1=u$). In other words, $\left(  v,v-1,v-2,\ldots,u\right)
=\left(  i_{1},i_{2},\ldots,i_{k}\right)  $.

We have shown that $v,v-1,v-2,\ldots,u$ are $v-u+1$ distinct elements of the
set $\left\{  1,2,\ldots,n\right\}  $. In other words, $i_{1},i_{2}%
,\ldots,i_{k}$ are $k$ distinct elements of the set $\left[  n\right]  $
(since $\left(  v,v-1,v-2,\ldots,u\right)  =\left(  i_{1},i_{2},\ldots
,i_{k}\right)  $ and $v-u+1=k$ and $\left\{  1,2,\ldots,n\right\}  =\left[
n\right]  $). Hence, Exercise \ref{exe.perm.c=ttt} yields%
\begin{equation}
\operatorname*{cyc}\nolimits_{i_{1},i_{2},\ldots,i_{k}}=t_{i_{1},i_{2}}\circ
t_{i_{2},i_{3}}\circ\cdots\circ t_{i_{k-1},i_{k}}.
\label{pf.cor.sol.perm.c=ttt.cuv1.b.2}%
\end{equation}

For each $j\in\left\{  1,2,\ldots,k-1\right\}  $, we have
\begin{equation}
t_{i_{j},i_{j+1}}=s_{v-j}. \label{pf.cor.sol.perm.c=ttt.cuv1.b.3}%
\end{equation}

[\textit{Proof of (\ref{pf.cor.sol.perm.c=ttt.cuv1.b.3}):} Let $j\in\left\{
1,2,\ldots,k-1\right\}  $. Thus, $j\in\left\{  1,2,\ldots,k-1\right\}
\subseteq\left\{  1,2,\ldots,k\right\}  $. Hence, the definition of $i_{j}$
yields $i_{j}=v-j+1$. Moreover, from $j\in\left\{  1,2,\ldots,k-1\right\}  $,
we obtain $j+1\in\left\{  2,3,\ldots,k\right\}  \subseteq\left\{
1,2,\ldots,k\right\}  $. Hence, the definition of $i_{j+1}$ yields
$i_{j+1}=v-\left(  j+1\right)  +1=v-j$. But from $j\in\left\{  1,2,\ldots
,k-1\right\}  $, we obtain%
\[
j\leq\underbrace{k}_{=v-u+1}-1=v-u+1-1=v-\underbrace{u}_{\geq1}\leq v-1.
\]
Thus, $v-j\geq1$. Combining this with%
\[
\underbrace{v}_{\leq n}-\underbrace{j}_{\substack{\geq1\\\text{(since }%
j\in\left\{  1,2,\ldots,k-1\right\}  \text{)}}}\leq n-1,
\]
we obtain $v-j\in\left\{  1,2,\ldots,n-1\right\}  $. Hence, Lemma
\ref{lem.sol.perm.c=ttt.s=t} (applied to $i=v-j$) yields $t_{v-j+1,v-j}%
=s_{v-j}$.

Now, from $i_{j}=v-j+1$ and $i_{j+1}=v-j$, we obtain $t_{i_{j},i_{j+1}%
}=t_{v-j+1,v-j}=s_{v-j}$. This proves (\ref{pf.cor.sol.perm.c=ttt.cuv1.b.3}).]

Now, recall that $\left(  i_{1},i_{2},\ldots,i_{k}\right)  =\left(
v,v-1,v-2,\ldots,u\right)  $. Hence,%
\[
\operatorname*{cyc}\nolimits_{i_{1},i_{2},\ldots,i_{k}}=\operatorname*{cyc}%
\nolimits_{v,v-1,v-2,\ldots,u}=c_{u,v}%
\]
(since $c_{u,v}=\operatorname*{cyc}\nolimits_{v,v-1,v-2,\ldots,u}$ (by the
definition of $c_{u,v}$)). Comparing this with%
\begin{align*}
\operatorname*{cyc}\nolimits_{i_{1},i_{2},\ldots,i_{k}}  &  =t_{i_{1},i_{2}%
}\circ t_{i_{2},i_{3}}\circ\cdots\circ t_{i_{k-1},i_{k}}%
\ \ \ \ \ \ \ \ \ \ \left(  \text{by (\ref{pf.cor.sol.perm.c=ttt.cuv1.b.2}%
)}\right) \\
&  =s_{v-1}\circ s_{v-2}\circ\cdots\circ s_{v-\left(  k-1\right)  }\\
&  \ \ \ \ \ \ \ \ \ \ \left(  \text{since each }j\in\left\{  1,2,\ldots
,k-1\right\}  \text{ satisfies }t_{i_{j},i_{j+1}}=s_{v-j}\text{ (by
(\ref{pf.cor.sol.perm.c=ttt.cuv1.b.3}))}\right) \\
&  =s_{v-1}\circ s_{v-2}\circ\cdots\circ s_{u}\ \ \ \ \ \ \ \ \ \ \left(
\text{since }v-\left(  k-1\right)  =v-k+1=u\right)  ,
\end{align*}
we obtain $c_{u,v}=s_{v-1}\circ s_{v-2}\circ\cdots\circ s_{u}$. This proves
Corollary \ref{cor.sol.perm.c=ttt.cuv1} \textbf{(b)}.
\end{verlong}
\end{proof}

\subsection{Solution to Exercise \ref{exe.perm.cycles}}

\begin{proof}
[Solution to Exercise \ref{exe.perm.cycles}.]\textbf{(a)} This proof is going
to be long, but most of it will be spent unraveling the notations. If you find
Exercise \ref{exe.perm.cycles} \textbf{(a)} obvious, don't let this proof cast
doubt on your understanding.

Let $\sigma\in S_{n}$. Let $i_{1},i_{2},\ldots,i_{k}$ be $k$ distinct elements
of $\left[  n\right]  $.

\begin{vershort}
The map $\sigma$ is a permutation (since $\sigma\in S_{n}$), and therefore
bijective. Hence, in particular, $\sigma$ is injective.

For every $p\in\left\{  1,2,\ldots,k\right\}  $, let $j_{p}$ be the element
$\sigma\left(  i_{p}\right)  \in\left[  n\right]  $. Then, $\left(
j_{1},j_{2},\ldots,j_{k}\right)  =\left(  \sigma\left(  i_{1}\right)
,\sigma\left(  i_{2}\right)  ,\ldots,\sigma\left(  i_{k}\right)  \right)  $.
Furthermore, $j_{1},j_{2},\ldots,j_{k}$ are $k$ distinct elements of $\left[
n\right]  $\ \ \ \ \footnote{\textit{Proof.} The map $\sigma$ is injective.
Hence, the elements $\sigma\left(  i_{1}\right)  ,\sigma\left(  i_{2}\right)
,\ldots,\sigma\left(  i_{k}\right)  $ are distinct (since the elements
$i_{1},i_{2},\ldots,i_{k}$ are distinct). Thus, $\sigma\left(  i_{1}\right)
,\sigma\left(  i_{2}\right)  ,\ldots,\sigma\left(  i_{k}\right)  $ are $k$
distinct elements of $\left[  n\right]  $. In other words, $j_{1},j_{2}%
,\ldots,j_{k}$ are $k$ distinct elements of $\left[  n\right]  $ (since
$\left(  j_{1},j_{2},\ldots,j_{k}\right)  =\left(  \sigma\left(  i_{1}\right)
,\sigma\left(  i_{2}\right)  ,\ldots,\sigma\left(  i_{k}\right)  \right)  $%
).}. Therefore, $\operatorname*{cyc}\nolimits_{j_{1},j_{2},\ldots,j_{k}}$ is a
well-defined permutation in $S_{n}$.

We have defined $\operatorname*{cyc}\nolimits_{i_{1},i_{2},\ldots,i_{k}}$ to
be the permutation in $S_{n}$ which sends $i_{1},i_{2},\ldots,i_{k}$ to
$i_{2},i_{3},\ldots,i_{k},i_{1}$, respectively, while leaving all other
elements of $\left[  n\right]  $ fixed. Therefore:

\begin{itemize}
\item The permutation $\operatorname*{cyc}\nolimits_{i_{1},i_{2},\ldots,i_{k}%
}$ sends $i_{1},i_{2},\ldots,i_{k}$ to $i_{2},i_{3},\ldots,i_{k},i_{1}$,
respectively. In other words,%
\begin{equation}
\operatorname*{cyc}\nolimits_{i_{1},i_{2},\ldots,i_{k}}\left(  i_{p}\right)
=i_{p+1}\ \ \ \ \ \ \ \ \ \ \text{for every }p\in\left\{  1,2,\ldots
,k\right\}  , \label{sol.perm.cycles.short.a.cyc-manifest.1}%
\end{equation}
where $i_{k+1}$ means $i_{1}$.

\item The permutation $\operatorname*{cyc}\nolimits_{i_{1},i_{2},\ldots,i_{k}%
}$ leaves all other elements of $\left[  n\right]  $ fixed (where
\textquotedblleft other\textquotedblright\ means \textquotedblleft other than
$i_{1},i_{2},\ldots,i_{k}$\textquotedblright). In other words,%
\begin{equation}
\operatorname*{cyc}\nolimits_{i_{1},i_{2},\ldots,i_{k}}\left(  q\right)
=q\ \ \ \ \ \ \ \ \ \ \text{for every }q\in\left[  n\right]  \setminus\left\{
i_{1},i_{2},\ldots,i_{k}\right\}  .
\label{sol.perm.cycles.short.a.cyc-manifest.2}%
\end{equation}

\end{itemize}

Similarly, we can say the same about $\operatorname*{cyc}\nolimits_{j_{1}%
,j_{2},\ldots,j_{k}}$:

\begin{itemize}
\item We have%
\begin{equation}
\operatorname*{cyc}\nolimits_{j_{1},j_{2},\ldots,j_{k}}\left(  j_{p}\right)
=j_{p+1}\ \ \ \ \ \ \ \ \ \ \text{for every }p\in\left\{  1,2,\ldots
,k\right\}  , \label{sol.perm.cycles.short.a.cyc-manifest.1'}%
\end{equation}
where $j_{k+1}$ means $j_{1}$.

\item We have%
\begin{equation}
\operatorname*{cyc}\nolimits_{j_{1},j_{2},\ldots,j_{k}}\left(  q\right)
=q\ \ \ \ \ \ \ \ \ \ \text{for every }q\in\left[  n\right]  \setminus\left\{
j_{1},j_{2},\ldots,j_{k}\right\}  .
\label{sol.perm.cycles.short.a.cyc-manifest.2'}%
\end{equation}

\end{itemize}

In the following, we shall use the notation $i_{k+1}$ as a synonym for $i_{1}%
$, and the notation $j_{k+1}$ as a synonym for $j_{1}$. Then,%
\begin{equation}
j_{p}=\sigma\left(  i_{p}\right)  \ \ \ \ \ \ \ \ \ \ \text{for every }%
p\in\left\{  1,2,\ldots,k+1\right\}  . \label{sol.perm.cycles.short.a.jp}%
\end{equation}
(Indeed, in the case when $p\in\left\{  1,2,\ldots,k\right\}  $, this follows
from the definition of $j_{p}$; but in the remaining case when $p=k+1$, it
follows from $j_{k+1}=j_{1}=\sigma\left(  \underbrace{i_{1}}_{=i_{k+1}%
}\right)  =\sigma\left(  i_{k+1}\right)  $.)

Now, let us show that
\begin{equation}
\left(  \sigma\circ\operatorname*{cyc}\nolimits_{i_{1},i_{2},\ldots,i_{k}%
}\circ\sigma^{-1}\right)  \left(  q\right)  =\operatorname*{cyc}%
\nolimits_{j_{1},j_{2},\ldots,j_{k}}\left(  q\right)
\label{sol.perm.cycles.short.a.qq}%
\end{equation}
for every $q\in\left[  n\right]  $.

[\textit{Proof of (\ref{sol.perm.cycles.short.a.qq}):} Let $q\in\left[
n\right]  $. We must prove (\ref{sol.perm.cycles.short.a.qq}). We are in one
of the following two cases:

\textit{Case 1:} We have $q\in\left\{  j_{1},j_{2},\ldots,j_{k}\right\}  $.

\textit{Case 2:} We have $q\notin\left\{  j_{1},j_{2},\ldots,j_{k}\right\}  $.

Let us first consider Case 1. In this case, we have $q\in\left\{  j_{1}%
,j_{2},\ldots,j_{k}\right\}  $. Thus, $q=j_{p}$ for some $p\in\left\{
1,2,\ldots,k\right\}  $. Consider this $p$. Clearly, $p+1\in\left\{
2,3,\ldots,k+1\right\}  \subseteq\left\{  1,2,\ldots,k+1\right\}  $. Hence,
applying (\ref{sol.perm.cycles.short.a.jp}) to $p+1$ instead of $p$, we obtain
$j_{p+1}=\sigma\left(  i_{p+1}\right)  $. But $q=j_{p}=\sigma\left(
i_{p}\right)  $ (by the definition of $j_{p}$) and thus $\sigma^{-1}\left(
q\right)  =i_{p}$. Hence,
\begin{align*}
\left(  \sigma\circ\operatorname*{cyc}\nolimits_{i_{1},i_{2},\ldots,i_{k}%
}\circ\sigma^{-1}\right)  \left(  q\right)   &  =\sigma\left(
\operatorname*{cyc}\nolimits_{i_{1},i_{2},\ldots,i_{k}}\left(
\underbrace{\sigma^{-1}\left(  q\right)  }_{=i_{p}}\right)  \right)
=\sigma\left(  \underbrace{\operatorname*{cyc}\nolimits_{i_{1},i_{2}%
,\ldots,i_{k}}\left(  i_{p}\right)  }_{\substack{=i_{p+1}\\\text{(by
(\ref{sol.perm.cycles.short.a.cyc-manifest.1}))}}}\right) \\
&  =\sigma\left(  i_{p+1}\right)  .
\end{align*}
Compared with%
\begin{align*}
\operatorname*{cyc}\nolimits_{j_{1},j_{2},\ldots,j_{k}}\left(  \underbrace{q}%
_{=j_{p}}\right)   &  =\operatorname*{cyc}\nolimits_{j_{1},j_{2},\ldots,j_{k}%
}\left(  j_{p}\right)  =j_{p+1}\ \ \ \ \ \ \ \ \ \ \left(  \text{by
(\ref{sol.perm.cycles.short.a.cyc-manifest.1'})}\right) \\
&  =\sigma\left(  i_{p+1}\right)  ,
\end{align*}
this yields $\left(  \sigma\circ\operatorname*{cyc}\nolimits_{i_{1}%
,i_{2},\ldots,i_{k}}\circ\sigma^{-1}\right)  \left(  q\right)
=\operatorname*{cyc}\nolimits_{j_{1},j_{2},\ldots,j_{k}}\left(  q\right)  $.
Thus, (\ref{sol.perm.cycles.short.a.qq}) is proven in Case 1.

Let us now consider Case 2. In this case, we have $q\notin\left\{  j_{1}%
,j_{2},\ldots,j_{k}\right\}  $. Hence, $\sigma^{-1}\left(  q\right)
\notin\left\{  i_{1},i_{2},\ldots,i_{k}\right\}  $%
\ \ \ \ \footnote{\textit{Proof.} Assume the contrary. Thus, $\sigma
^{-1}\left(  q\right)  \in\left\{  i_{1},i_{2},\ldots,i_{k}\right\}  $. In
other words, there exists a $p\in\left\{  1,2,\ldots,k\right\}  $ such that
$\sigma^{-1}\left(  q\right)  =i_{p}$. Consider this $p$. We have $\sigma
^{-1}\left(  q\right)  =i_{p}$, thus $q=\sigma\left(  i_{p}\right)  =j_{p}$
(since $j_{p}$ is defined as $\sigma\left(  i_{p}\right)  $). Thus,
$q=j_{p}\in\left\{  j_{1},j_{2},\ldots,j_{k}\right\}  $; but this contradicts
$q\notin\left\{  j_{1},j_{2},\ldots,j_{k}\right\}  $. This contradiction shows
that our assumption was wrong, qed.}. Thus, $\sigma^{-1}\left(  q\right)
\in\left[  n\right]  \setminus\left\{  i_{1},i_{2},\ldots,i_{k}\right\}  $.
Therefore, (\ref{sol.perm.cycles.short.a.cyc-manifest.2}) (applied to
$\sigma^{-1}\left(  q\right)  $ instead of $q$) yields $\operatorname*{cyc}%
\nolimits_{i_{1},i_{2},\ldots,i_{k}}\left(  \sigma^{-1}\left(  q\right)
\right)  =\sigma^{-1}\left(  q\right)  $. Hence,%
\begin{equation}
\left(  \sigma\circ\operatorname*{cyc}\nolimits_{i_{1},i_{2},\ldots,i_{k}%
}\circ\sigma^{-1}\right)  \left(  q\right)  =\sigma\left(
\underbrace{\operatorname*{cyc}\nolimits_{i_{1},i_{2},\ldots,i_{k}}\left(
\sigma^{-1}\left(  q\right)  \right)  }_{=\sigma^{-1}\left(  q\right)
}\right)  =\sigma\left(  \sigma^{-1}\left(  q\right)  \right)  =q.
\label{sol.perm.cycles.short.a.qq.pf.5}%
\end{equation}
On the other hand, $q\in\left[  n\right]  \setminus\left\{  j_{1},j_{2}%
,\ldots,j_{k}\right\}  $ (since $q\notin\left\{  j_{1},j_{2},\ldots
,j_{k}\right\}  $) and thus $\operatorname*{cyc}\nolimits_{j_{1},j_{2}%
,\ldots,j_{k}}\left(  q\right)  =q$ (by
(\ref{sol.perm.cycles.short.a.cyc-manifest.2'})). Compared with
(\ref{sol.perm.cycles.short.a.qq.pf.5}), this yields \newline$\left(
\sigma\circ\operatorname*{cyc}\nolimits_{i_{1},i_{2},\ldots,i_{k}}\circ
\sigma^{-1}\right)  \left(  q\right)  =\operatorname*{cyc}\nolimits_{j_{1}%
,j_{2},\ldots,j_{k}}\left(  q\right)  $. Thus,
(\ref{sol.perm.cycles.short.a.qq}) is proven in Case 2.

We have now proven (\ref{sol.perm.cycles.short.a.qq}) in each of the two Cases
1 and 2. Therefore, (\ref{sol.perm.cycles.short.a.qq}) always holds.]

From this, we conclude that $\sigma\circ\operatorname*{cyc}\nolimits_{i_{1}%
,i_{2},\ldots,i_{k}}\circ\sigma^{-1}=\operatorname*{cyc}\nolimits_{j_{1}%
,j_{2},\ldots,j_{k}}$ (because both $\sigma\circ\operatorname*{cyc}%
\nolimits_{i_{1},i_{2},\ldots,i_{k}}\circ\sigma^{-1}$ and $\operatorname*{cyc}%
\nolimits_{j_{1},j_{2},\ldots,j_{k}}$ are maps from $\left[  n\right]  $ to
$\left[  n\right]  $). Hence,%
\[
\sigma\circ\operatorname*{cyc}\nolimits_{i_{1},i_{2},\ldots,i_{k}}\circ
\sigma^{-1}=\operatorname*{cyc}\nolimits_{j_{1},j_{2},\ldots,j_{k}%
}=\operatorname*{cyc}\nolimits_{\sigma\left(  i_{1}\right)  ,\sigma\left(
i_{2}\right)  ,\ldots,\sigma\left(  i_{k}\right)  }%
\]
(since $\left(  j_{1},j_{2},\ldots,j_{k}\right)  =\left(  \sigma\left(
i_{1}\right)  ,\sigma\left(  i_{2}\right)  ,\ldots,\sigma\left(  i_{k}\right)
\right)  $). This solves Exercise \ref{exe.perm.cycles} \textbf{(a)}.

\textbf{(b)} This is again a straightforward argument whose complexity stems
only from the number of cases that need to be considered. We shall try to
reduce the amount of brainless verification using some tricks, although at the
cost of making parts of the solution appear unmotivated.

Let $p\in\left\{  0,1,\ldots,n-k\right\}  $. Then, the elements
$p+1,p+2,\ldots,p+k$ belong to $\left\{  1,2,\ldots,n\right\}  =\left[
n\right]  $. Hence, $\operatorname*{cyc}\nolimits_{p+1,p+2,\ldots,p+k}$ is
well-defined. We let $\sigma$ denote the permutation $\operatorname*{cyc}%
\nolimits_{p+1,p+2,\ldots,p+k}$.

The permutation $\sigma=\operatorname*{cyc}\nolimits_{p+1,p+2,\ldots,p+k}$ is
defined to be the permutation in $S_{n}$ which sends $p+1,p+2,\ldots,p+k$ to
$p+2,p+3,\ldots,p+k,p+1$, respectively, while leaving all other elements of
$\left[  n\right]  $ fixed. Therefore:

\begin{itemize}
\item The permutation $\sigma$ sends $p+1,p+2,\ldots,p+k$ to $p+2,p+3,\ldots
,p+k,p+1$, respectively. In other words, we have%
\begin{equation}
\left(  \sigma\left(  p+i\right)  =p+\left(  i+1\right)
\ \ \ \ \ \ \ \ \ \ \text{for every }i\in\left\{  1,2,\ldots,k-1\right\}
\right)  \label{sol.perm.cycles.short.b.sigma-manifest.1a}%
\end{equation}
and%
\begin{equation}
\sigma\left(  p+k\right)  =p+1.
\label{sol.perm.cycles.short.b.sigma-manifest.1b}%
\end{equation}

\item The permutation $\sigma$ leaves all other elements of $\left[  n\right]
$ fixed (where \textquotedblleft other\textquotedblright\ means
\textquotedblleft other than $p+1,p+2,\ldots,p+k$\textquotedblright). In other
words,
\begin{equation}
\sigma\left(  q\right)  =q\ \ \ \ \ \ \ \ \ \ \text{for every }q\in\left[
n\right]  \setminus\left\{  p+1,p+2,\ldots,p+k\right\}  .
\label{sol.perm.cycles.short.b.sigma-manifest.2}%
\end{equation}

\end{itemize}

We can now observe that
\begin{equation}
\sigma\left(  q\right)  \geq q\ \ \ \ \ \ \ \ \ \ \text{for every }q\in\left[
n\right]  \text{ satisfying }q\neq p+k \label{sol.perm.cycles.short.b.sigma.1}%
\end{equation}
\footnote{\textit{Proof of (\ref{sol.perm.cycles.short.b.sigma.1}):} Let
$q\in\left[  n\right]  $ be such that $q\neq p+k$. We must prove that
$\sigma\left(  q\right)  \geq q$.
\par
We are in one of the following two cases:
\par
\textit{Case 1:} We have $q\in\left\{  p+1,p+2,\ldots,p+k\right\}  $.
\par
\textit{Case 2:} We have $q\notin\left\{  p+1,p+2,\ldots,p+k\right\}  $.
\par
Let us first consider Case 1. In this case, we have $q\in\left\{
p+1,p+2,\ldots,p+k\right\}  $. Hence, $q=p+i$ for some $i\in\left\{
1,2,\ldots,k\right\}  $. Consider this $i$. We have $q\neq p+k$, so that
$q-p\neq k$ and thus $k\neq q-p=i$ (since $q=p+i$). Therefore, $i\neq k$.
Combined with $i\in\left\{  1,2,\ldots,k\right\}  $, this shows that
$i\in\left\{  1,2,\ldots,k\right\}  \setminus\left\{  k\right\}  =\left\{
1,2,\ldots,k-1\right\}  $. Hence,
(\ref{sol.perm.cycles.short.b.sigma-manifest.1a}) yields $\sigma\left(
p+i\right)  =p+\underbrace{\left(  i+1\right)  }_{\geq i}\geq p+i=q$. Hence,
$\sigma\left(  \underbrace{q}_{=p+i}\right)  =\sigma\left(  p+i\right)  \geq
q$. We thus have proven $\sigma\left(  q\right)  \geq q$ in Case 1.
\par
Let us now consider Case 2. In this case, we have $q\notin\left\{
p+1,p+2,\ldots,p+k\right\}  $. Hence, $q\in\left[  n\right]  \setminus\left\{
p+1,p+2,\ldots,p+k\right\}  $. Therefore, $\sigma\left(  q\right)  =q$ (by
(\ref{sol.perm.cycles.short.b.sigma-manifest.2})). Thus, $\sigma\left(
q\right)  \geq q$ is proven in Case 2.
\par
We have now proven $\sigma\left(  q\right)  \geq q$ in both Cases 1 and 2.
Therefore, $\sigma\left(  q\right)  \geq q$ always holds. This proves
(\ref{sol.perm.cycles.short.b.sigma.1}).}. Furthermore,%
\begin{equation}
\sigma\left(  q\right)  \leq q+1\ \ \ \ \ \ \ \ \ \ \text{for every }%
q\in\left[  n\right]  \label{sol.perm.cycles.short.b.sigma.2}%
\end{equation}
\footnote{\textit{Proof of (\ref{sol.perm.cycles.short.b.sigma.2}):} Let
$q\in\left[  n\right]  $. We must prove that $\sigma\left(  q\right)  \leq
q+1$.
\par
Assume the contrary (for the sake of contradiction). Thus, $\sigma\left(
q\right)  >q+1$. If we had $q\in\left[  n\right]  \setminus\left\{
p+1,p+2,\ldots,p+k\right\}  $, then we would have $\sigma\left(  q\right)  =q$
(by (\ref{sol.perm.cycles.short.b.sigma-manifest.2})), which would contradict
$\sigma\left(  q\right)  >q+1>q$. Thus, we cannot have $q\in\left[  n\right]
\setminus\left\{  p+1,p+2,\ldots,p+k\right\}  $. We therefore have%
\[
q\in\left[  n\right]  \setminus\left(  \left[  n\right]  \setminus\left\{
p+1,p+2,\ldots,p+k\right\}  \right)  \subseteq\left\{  p+1,p+2,\ldots
,p+k\right\}  .
\]
Hence, $q=p+i$ for some $i\in\left\{  1,2,\ldots,k\right\}  $. Consider this
$i$. Clearly, $i\geq1$ and $i\leq k$.
\par
If we had $i=k$, then we would have%
\begin{align*}
\sigma\left(  \underbrace{q}_{=p+i}\right)   &  =\sigma\left(
p+\underbrace{i}_{=k}\right)  =\sigma\left(  p+k\right)  =p+\underbrace{1}%
_{\leq i}\ \ \ \ \ \ \ \ \ \ \left(  \text{by
(\ref{sol.perm.cycles.short.b.sigma-manifest.1b})}\right) \\
&  \leq p+i=q,
\end{align*}
which would contradict $\sigma\left(  q\right)  >q+1>q$. Thus, we cannot have
$i=k$. Hence, $i\in\left\{  1,2,\ldots,k-1\right\}  $ (since $i\geq1$ and
$i\leq k$). Therefore, (\ref{sol.perm.cycles.short.b.sigma-manifest.1a})
yields $\sigma\left(  p+i\right)  =p+\left(  i+1\right)  =\underbrace{p+i}%
_{=q}+1=q+1$. This contradicts $\sigma\left(  \underbrace{p+i}_{=q}\right)
=\sigma\left(  q\right)  >q+1$. This contradiction proves that our assumption
was wrong. Hence, $\sigma\left(  q\right)  \leq q+1$ is proven.}.

Now, set%
\[
A=\left\{  \left(  p+h,p+k\right)  \ \mid\ h\in\left\{  1,2,\ldots
,k-1\right\}  \right\}  .
\]
In other words,%
\[
A=\left\{  \left(  p+1,p+k\right)  ,\left(  p+2,p+k\right)  ,\ldots,\left(
p+\left(  k-1\right)  ,p+k\right)  \right\}  .
\]
Thus, the set $A$ has $k-1$ elements (since the $k-1$ pairs \newline$\left(
p+1,p+k\right)  ,\left(  p+2,p+k\right)  ,\ldots,\left(  p+\left(  k-1\right)
,p+k\right)  $ are clearly distinct). In other words, $\left\vert A\right\vert
=k-1$.

Now, let $\operatorname*{Inv}\sigma$ denote the set of all inversions of
$\sigma$. Recall that $\ell\left(  \sigma\right)  $ was defined as the number
of inversions of $\sigma$. In other words, $\ell\left(  \sigma\right)
=\left\vert \operatorname*{Inv}\sigma\right\vert $.

But $A\subseteq\operatorname*{Inv}\sigma$\ \ \ \ \footnote{\textit{Proof.} Let
$c\in A$. We shall show that $c\in\operatorname*{Inv}\sigma$.
\par
We have $c\in A=\left\{  \left(  p+h,p+k\right)  \ \mid\ h\in\left\{
1,2,\ldots,k-1\right\}  \right\}  $. Hence, $c$ can be written in the form
$c=\left(  p+h,p+k\right)  $ for some $h\in\left\{  1,2,\ldots,k-1\right\}  $.
Consider this $h$. From $h\in\left\{  1,2,\ldots,k-1\right\}  $, we obtain
$1\leq h\leq k-1$, so that $h\leq k-1<k$. Thus, $p+\underbrace{h}_{<k}<p+k$.
Moreover, $1\leq p+h$ (since $\underbrace{p}_{\geq0}+\underbrace{h}_{\geq
1}\geq0+1=1$) and $p+k\leq n$ (since $p\leq n-k$). Thus, $1\leq p+h<p+k\leq
n$.
\par
Applying (\ref{sol.perm.cycles.short.b.sigma-manifest.1a}) to $i=h$, we obtain
$\sigma\left(  p+h\right)  =p+\left(  \underbrace{h}_{\geq1>0}+1\right)
>p+\left(  0+1\right)  =p+1=\sigma\left(  p+k\right)  $ (by
(\ref{sol.perm.cycles.short.b.sigma-manifest.1b})).
\par
Now, $\left(  p+h,p+k\right)  $ is a pair of integers satisfying $1\leq
p+h<p+k\leq n$ and $\sigma\left(  p+h\right)  >\sigma\left(  p+k\right)  $. In
other words, $\left(  p+h,p+k\right)  $ is a pair of integers $\left(
i,j\right)  $ satisfying $1\leq i<j\leq n$ and $\sigma\left(  i\right)
>\sigma\left(  j\right)  $. In other words, $\left(  p+h,p+k\right)  $ is an
inversion of $\sigma$ (by the definition of \textquotedblleft inversion of
$\sigma$\textquotedblright). In other words, $\left(  p+h,p+k\right)
\in\operatorname*{Inv}\sigma$. Thus, $c=\left(  p+h,p+k\right)  \in
\operatorname*{Inv}\sigma$.
\par
Now, let us forget that we fixed $c$. We thus have shown that $c\in
\operatorname*{Inv}\sigma$ for every $c\in A$. In other words, $A\subseteq
\operatorname*{Inv}\sigma$, qed.} and $\operatorname*{Inv}\sigma\subseteq
A$\ \ \ \ \footnote{\textit{Proof.} Let $c\in\operatorname*{Inv}\sigma$. We
shall show that $c\in A$.
\par
We have $c\in\operatorname*{Inv}\sigma$. In other words, $c$ is an inversion
of $\sigma$. In other words, $c$ is a pair $\left(  i,j\right)  $ of integers
satisfying $1\leq i<j\leq n$ and $\sigma\left(  i\right)  >\sigma\left(
j\right)  $. Consider this $\left(  i,j\right)  $.
\par
We have $\sigma\left(  i\right)  >\sigma\left(  j\right)  $, so that
$\sigma\left(  i\right)  \geq\sigma\left(  j\right)  +1$ (since $\sigma\left(
i\right)  $ and $\sigma\left(  j\right)  $ are integers). Also, $i<j$, so that
$i\leq j-1$ (since $i$ and $j$ are integers). In other words, $i+1\leq j$. But
(\ref{sol.perm.cycles.short.b.sigma.2}) (applied to $q=i$) yields
$\sigma\left(  i\right)  \leq i+1\leq j$.
\par
Let us first show that $j=p+k$. Indeed, let us assume the contrary (for the
sake of contradiction). Thus, $j\neq p+k$. Hence, $\sigma\left(  j\right)
\geq j$ (by (\ref{sol.perm.cycles.short.b.sigma.1}), applied to $q=j$). Now,
$\sigma\left(  i\right)  >\sigma\left(  j\right)  \geq j$. This contradicts
$\sigma\left(  i\right)  \leq j$. This contradiction shows that our assumption
was wrong. Hence, $j=p+k$ is proven.
\par
Now, $\sigma\left(  \underbrace{j}_{=p+k}\right)  =\sigma\left(  p+k\right)
=p+1$ (by (\ref{sol.perm.cycles.short.b.sigma-manifest.1b})). Hence,
$\sigma\left(  i\right)  >\sigma\left(  j\right)  =p+1$. Therefore,
$p+1<\sigma\left(  i\right)  \leq i+1$. Subtracting $1$ from both sides of
this inequality, we obtain $p<i$. Hence, $i>p$, so that $i\geq p+1$ (since $i$
and $p$ are integers). Combined with $i<j=p+k$, this yields $i\in\left\{
p+1,p+2,\ldots,p+k-1\right\}  $. Thus, $i-p\in\left\{  1,2,\ldots,k-1\right\}
$.
\par
So we know that the element $i-p\in\left\{  1,2,\ldots,k-1\right\}  $
satisfies $c=\left(  \underbrace{i}_{=p+\left(  i-p\right)  },\underbrace{j}%
_{=p+k}\right)  =\left(  p+\left(  i-p\right)  ,p+k\right)  $. Hence, there
exists an $h\in\left\{  1,2,\ldots,k-1\right\}  $ such that $c=\left(
p+h,p+k\right)  $ (namely, $h=i-p$). Thus,%
\[
c\in\left\{  \left(  p+h,p+k\right)  \ \mid\ h\in\left\{  1,2,\ldots
,k-1\right\}  \right\}  =A.
\]
\par
Now, let us forget that we fixed $c$. We thus have proven that $c\in A$ for
every $c\in\operatorname*{Inv}\sigma$. In other words, $\operatorname*{Inv}%
\sigma\subseteq A$, qed.}. Combining these two relations, we obtain
$A=\operatorname*{Inv}\sigma$. Hence, $\left\vert A\right\vert =\left\vert
\operatorname*{Inv}\sigma\right\vert $. Compared with $\ell\left(
\sigma\right)  =\left\vert \operatorname*{Inv}\sigma\right\vert $, this yields
$\ell\left(  \sigma\right)  =\left\vert A\right\vert =k-1$. This rewrites as
$\ell\left(  \operatorname*{cyc}\nolimits_{p+1,p+2,\ldots,p+k}\right)  =k-1$
(since $\sigma=\operatorname*{cyc}\nolimits_{p+1,p+2,\ldots,p+k}$). This
solves Exercise \ref{exe.perm.cycles} \textbf{(b)}.

\textbf{(c)} This one is tricky. Let $i_{1},i_{2},\ldots,i_{k}$ be $k$
distinct elements of $\left[  n\right]  $. We extend the $k$-tuple $\left(
i_{1},i_{2},\ldots,i_{k}\right)  $ to an infinite sequence $\left(
i_{1},i_{2},i_{3},\ldots\right)  $ of elements of $\left[  n\right]  $ by
setting%
\[
\left(  i_{u}=i_{\left(  \text{the element }u^{\prime}\in\left\{
1,2,\ldots,k\right\}  \text{ satisfying }u^{\prime}\equiv u\operatorname{mod}%
k\right)  }\ \ \ \ \ \ \ \ \ \ \text{for every }u\geq1\right)  .
\]
This sequence $\left(  i_{1},i_{2},i_{3},\ldots\right)  $ is periodic with
period $k$. In other words,%
\begin{equation}
i_{u}=i_{u+k}\ \ \ \ \ \ \ \ \ \ \text{for every }u\geq1.
\label{sol.perm.cycles.short.c.period}%
\end{equation}
From this, it is easy to obtain that%
\[
i_{1}+i_{2}+\cdots+i_{k}=i_{2}+i_{3}+\cdots+i_{k+1}=i_{3}+i_{4}+\cdots
+i_{k+2}=\cdots
\]
\footnote{\textit{Proof.} Every $u\in\left\{  1,2,3,\ldots\right\}  $
satisfies%
\begin{align*}
i_{u}+i_{u+1}+\cdots+i_{u+k-1}  &  =\underbrace{i_{u}}_{\substack{=i_{u+k}%
\\\text{(by (\ref{sol.perm.cycles.short.c.period}))}}}+\left(  i_{u+1}%
+i_{u+2}+\cdots+i_{u+k-1}\right) \\
&  =i_{u+k}+\left(  i_{u+1}+i_{u+2}+\cdots+i_{u+k-1}\right)  =\left(
i_{u+1}+i_{u+2}+\cdots+i_{u+k-1}\right)  +i_{u+k}\\
&  =i_{u+1}+i_{u+2}+\cdots+i_{u+k}.
\end{align*}
Thus, $i_{1}+i_{2}+\cdots+i_{k}=i_{2}+i_{3}+\cdots+i_{k+1}=i_{3}+i_{4}%
+\cdots+i_{k+2}=\cdots$, qed.}. Thus,%
\begin{equation}
i_{r+1}+i_{r+2}+\cdots+i_{r+k}=i_{1}+i_{2}+\cdots+i_{k}%
\ \ \ \ \ \ \ \ \ \ \text{for every }r\in\mathbb{N}.
\label{sol.perm.cycles.short.c.periodsum}%
\end{equation}

Let $\sigma=\operatorname*{cyc}\nolimits_{i_{1},i_{2},\ldots,i_{k}}$. Let
$\operatorname*{Inv}\sigma$ denote the set of all inversions of $\sigma$.
Then, $\ell\left(  \sigma\right)  =\left\vert \operatorname*{Inv}%
\sigma\right\vert $. (This can be seen as in the solution to Exercise
\ref{exe.perm.cycles} \textbf{(b)}.) Moreover, the definitions of the sequence
$\left(  i_{1},i_{2},i_{3},\ldots\right)  $ and of $\sigma$ show that%
\begin{equation}
\sigma\left(  i_{p}\right)  =i_{p+1}\ \ \ \ \ \ \ \ \ \ \text{for every }%
p\geq1. \label{sol.perm.cycles.short.c.sigma}%
\end{equation}

Now, fix $r\in\left\{  1,2,\ldots,k-1\right\}  $. We shall prove that
\begin{equation}
\text{there exists some }u\geq1\text{ such that }\left(  i_{u},i_{u+r}\right)
\in\operatorname*{Inv}\sigma. \label{sol.perm.cycles.short.c.mainclaim}%
\end{equation}

[\textit{Proof of (\ref{sol.perm.cycles.short.c.mainclaim}):} The sequence
$\left(  i_{1},i_{2},i_{3},\ldots\right)  $ is periodic with period $k$, and
its first $k$ entries $i_{1},i_{2},\ldots,i_{k}$ are distinct. Hence, each
entry of this sequence repeats itself each $k$ steps, but not more often.
Hence, every integer $u\geq1$ satisfies $i_{u+r}\neq i_{u}$ (since
$r\in\left\{  1,2,\ldots,k-1\right\}  $). In other words, every $u\geq1$
satisfies%
\begin{equation}
i_{u+r}-i_{u}\neq0. \label{sol.perm.cycles.short.c.mainclaim.pf.4}%
\end{equation}
The $k$-tuple $\left(  i_{1+r}-i_{1},i_{2+r}-i_{2},\ldots,i_{k+r}%
-i_{k}\right)  $ contains at least one positive entry\footnote{\textit{Proof.}
Assume the contrary. Thus, the $k$-tuple $\left(  i_{1+r}-i_{1},i_{2+r}%
-i_{2},\ldots,i_{k+r}-i_{k}\right)  $ contains no positive entries. In other
words, no $u\in\left\{  1,2,\ldots,k\right\}  $ satisfies $i_{u+r}-i_{u}>0$.
In other words, every $u\in\left\{  1,2,\ldots,k\right\}  $ satisfies
$i_{u+r}-i_{u}\leq0$. But (\ref{sol.perm.cycles.short.c.mainclaim.pf.4}) shows
that every $u\in\left\{  1,2,\ldots,k\right\}  $ satisfies $i_{u+r}-i_{u}%
\neq0$. Thus, every $u\in\left\{  1,2,\ldots,k\right\}  $ satisfies
$i_{u+r}-i_{u}<0$ (since $i_{u+r}-i_{u}\leq0$ and $i_{u+r}-i_{u}\neq0$).
Hence,%
\[
\sum_{u=1}^{k}\underbrace{\left(  i_{u+r}-i_{u}\right)  }_{<0}<\sum_{u=1}%
^{k}0=0.
\]
But this contradicts%
\begin{align*}
\sum_{u=1}^{k}\left(  \underbrace{i_{u+r}}_{=i_{r+u}}-i_{u}\right)   &
=\sum_{u=1}^{k}\left(  i_{r+u}-i_{u}\right)  =\left(  \sum_{u=1}^{k}%
i_{r+u}\right)  -\left(  \sum_{u=1}^{k}i_{u}\right) \\
&  =\left(  i_{r+1}+i_{r+2}+\cdots+i_{r+k}\right)  -\left(  i_{1}+i_{2}%
+\cdots+i_{k}\right)  =0\ \ \ \ \ \ \ \ \ \ \left(  \text{by
(\ref{sol.perm.cycles.short.c.periodsum})}\right)  .
\end{align*}
This contradiction shows that our assumption was wrong, qed.}, and at least
one negative entry\footnote{This is proven similarly.}. Hence, there exists at
least one $u\geq1$ such that $i_{u+r}-i_{u}>0$ but $i_{\left(  u+1\right)
+r}-i_{u+1}<0$\ \ \ \ \footnote{\textit{Proof.} Assume the contrary. Then,
there exists no $u\geq1$ such that $i_{u+r}-i_{u}>0$ but $i_{\left(
u+1\right)  +r}-i_{u+1}<0$. Hence,%
\begin{equation}
\text{every }u\geq1\text{ satisfying }i_{u+r}-i_{u}>0\text{ must satisfy
}i_{\left(  u+1\right)  +r}-i_{u+1}\geq0.
\label{sol.perm.cycles.short.c.mainclaim.pf.5o}%
\end{equation}
Therefore,%
\begin{equation}
\text{every }u\geq1\text{ satisfying }i_{u+r}-i_{u}>0\text{ must satisfy
}i_{\left(  u+1\right)  +r}-i_{u+1}>0
\label{sol.perm.cycles.short.c.mainclaim.pf.5}%
\end{equation}
(because (\ref{sol.perm.cycles.short.c.mainclaim.pf.5o}) shows that
$i_{\left(  u+1\right)  +r}-i_{u+1}\geq0$; but combining this with $i_{\left(
u+1\right)  +r}-i_{u+1}\neq0$ (which follows from
(\ref{sol.perm.cycles.short.c.mainclaim.pf.4}), applied to $u+1$ instead of
$u$), we obtain $i_{\left(  u+1\right)  +r}-i_{u+1}>0$).
\par
But there exists some $v\in\left\{  1,2,\ldots,k\right\}  $ satisfying
$i_{v+r}-i_{v}>0$ (since the $k$-tuple $\left(  i_{1+r}-i_{1},i_{2+r}%
-i_{2},\ldots,i_{k+r}-i_{k}\right)  $ contains at least one positive entry).
Consider this $v$. Then, we have $i_{v+r}-i_{v}>0$, therefore $i_{\left(
v+1\right)  +r}-i_{v+1}>0$ (by (\ref{sol.perm.cycles.short.c.mainclaim.pf.5}),
applied to $u=v$), therefore $i_{\left(  v+2\right)  +r}-i_{v+2}>0$ (by
(\ref{sol.perm.cycles.short.c.mainclaim.pf.5}), applied to $u=v+1$), therefore
$i_{\left(  v+3\right)  +r}-i_{v+3}>0$ (by
(\ref{sol.perm.cycles.short.c.mainclaim.pf.5}), applied to $u=v+2$), and so
on. Altogether, we thus obtain
\[
i_{h+r}-i_{h}>0\ \ \ \ \ \ \ \ \ \ \text{for every }h\geq v.
\]
In other words, $i_{h}<i_{h+r}$ for every $h\geq v$. Hence, $i_{v}%
<i_{v+r}<i_{v+2r}<i_{v+3r}<\cdots$. Thus, the numbers $i_{v},i_{v+r}%
,i_{v+2r},i_{v+3r},\ldots$ are pairwise distinct; hence, the sequence $\left(
i_{1},i_{2},i_{3},\ldots\right)  $ contains infinitely many distinct entries.
But this contradicts the fact that this sequence is periodic. This
contradiction proves that our assumption was wrong, qed.}. Consider this $u$.
We have $i_{u}<i_{u+r}$ (since $i_{u+r}-i_{u}>0$), so that $1\leq
i_{u}<i_{u+r}\leq n$. Also, (\ref{sol.perm.cycles.short.c.sigma}) (applied to
$p=u$) yields $\sigma\left(  i_{u}\right)  =i_{u+1}$. Moreover,
(\ref{sol.perm.cycles.short.c.sigma}) (applied to $p=u+r$) yields
$\sigma\left(  i_{u+r}\right)  =i_{u+r+1}=i_{\left(  u+1\right)  +r}$. Hence,
\begin{align*}
\sigma\left(  i_{u}\right)   &  =i_{u+1}>i_{\left(  u+1\right)  +r}%
\ \ \ \ \ \ \ \ \ \ \left(  \text{since }i_{\left(  u+1\right)  +r}%
-i_{u+1}<0\right) \\
&  =\sigma\left(  i_{u+r}\right)  .
\end{align*}
So we know that $\left(  i_{u},i_{u+r}\right)  $ is a pair of integers
satisfying $1\leq i_{u}<i_{u+r}\leq n$ and $\sigma\left(  i_{u}\right)
>\sigma\left(  i_{u+r}\right)  $. In other words, $\left(  i_{u}%
,i_{u+r}\right)  $ is an inversion of $\sigma$. In other words, $\left(
i_{u},i_{u+r}\right)  \in\operatorname*{Inv}\sigma$. Thus, we have found a
$u\geq1$ such that $\left(  i_{u},i_{u+r}\right)  \in\operatorname*{Inv}%
\sigma$. This proves (\ref{sol.perm.cycles.short.c.mainclaim}).]

Now, let us forget that we fixed $r$. We have shown that, for every
$r\in\left\{  1,2,\ldots,k-1\right\}  $, there exists some $u\geq1$ such that
$\left(  i_{u},i_{u+r}\right)  \in\operatorname*{Inv}\sigma$. Let us denote
this $u$ by $u_{r}$. Therefore, for every $r\in\left\{  1,2,\ldots
,k-1\right\}  $, we have found a $u_{r}\geq1$ such that $\left(  i_{u_{r}%
},i_{u_{r}+r}\right)  \in\operatorname*{Inv}\sigma$. The $k-1$ pairs%
\[
\left(  i_{u_{1}},i_{u_{1}+1}\right)  ,\ \left(  i_{u_{2}},i_{u_{2}+2}\right)
,\ \ldots,\ \left(  i_{u_{k-1}},i_{u_{k-1}+\left(  k-1\right)  }\right)
\]
are pairwise distinct\footnote{\textit{Proof.} Assume the contrary. Then,
there exist two distinct elements $x$ and $y$ of $\left\{  1,2,\ldots
,k-1\right\}  $ such that $\left(  i_{u_{x}},i_{u_{x}+x}\right)  =\left(
i_{u_{y}},i_{u_{y}+y}\right)  $. Consider these $x$ and $y$.
\par
We have $\left(  i_{u_{x}},i_{u_{x}+x}\right)  =\left(  i_{u_{y}},i_{u_{y}%
+y}\right)  $. In other words, $i_{u_{x}}=i_{u_{y}}$ and $i_{u_{x}+x}%
=i_{u_{y}+y}$. Since the numbers $i_{1},i_{2},\ldots,i_{k}$ are distinct (and
the sequence $\left(  i_{1},i_{2},i_{3},\ldots\right)  $ consists of these
numbers, repeated over and over), we obtain $u_{x}\equiv u_{y}%
\operatorname{mod}k$ from $i_{u_{x}}=i_{u_{y}}$, and we obtain $u_{x}+x\equiv
u_{y}+y\operatorname{mod}k$ from $i_{u_{x}+x}=i_{u_{y}+y}$. Subtracting the
congruence $u_{x}\equiv u_{y}\operatorname{mod}k$ from the congruence
$u_{x}+x\equiv u_{y}+y\operatorname{mod}k$, we obtain $x\equiv
y\operatorname{mod}k$. In light of $x,y\in\left\{  1,2,\ldots,k-1\right\}  $,
this shows that $x=y$. But this contradicts the fact that $x$ and $y$ are
distinct. This contradiction proves that our assumption was wrong, qed.}, and
all of them belong to $\operatorname*{Inv}\sigma$. Hence, the set
$\operatorname*{Inv}\sigma$ has at least $k-1$ elements. In other words,
$\left\vert \operatorname*{Inv}\sigma\right\vert \geq k-1$. Thus, $\ell\left(
\sigma\right)  =\left\vert \operatorname*{Inv}\sigma\right\vert \geq k-1$.
Since $\sigma=\operatorname*{cyc}\nolimits_{i_{1},i_{2},\ldots,i_{k}}$, this
rewrites as $\ell\left(  \operatorname*{cyc}\nolimits_{i_{1},i_{2}%
,\ldots,i_{k}}\right)  \geq k-1$. This solves Exercise \ref{exe.perm.cycles}
\textbf{(c)}.

\textbf{(d)} \textit{First solution to Exercise \ref{exe.perm.cycles}
\textbf{(d)}:} Let $i_{1},i_{2},\ldots,i_{k}$ be $k$ distinct elements of
$\left[  n\right]  $. Hence, Proposition \ref{prop.perms.lists} \textbf{(c)}
(applied to $\left(  p_{1},p_{2},\ldots,p_{k}\right)  =\left(  i_{1}%
,i_{2},\ldots,i_{k}\right)  $) yields that there exists a permutation
$\sigma\in S_{n}$ such that $\left(  i_{1},i_{2},\ldots,i_{k}\right)  =\left(
\sigma\left(  1\right)  ,\sigma\left(  2\right)  ,\ldots,\sigma\left(
k\right)  \right)  $. Consider such a $\sigma$.

Exercise \ref{exe.perm.cycles} \textbf{(b)} yields $\ell\left(
\operatorname*{cyc}\nolimits_{1,2,\ldots,k}\right)  =k-1$. But the definition
of $\left(  -1\right)  ^{\operatorname*{cyc}\nolimits_{1,2,\ldots,k}}$ yields
$\left(  -1\right)  ^{\operatorname*{cyc}\nolimits_{1,2,\ldots,k}}=\left(
-1\right)  ^{\ell\left(  \operatorname*{cyc}\nolimits_{1,2,\ldots,k}\right)
}=\left(  -1\right)  ^{k-1}$ (since $\ell\left(  \operatorname*{cyc}%
\nolimits_{1,2,\ldots,k}\right)  =k-1$).

Exercise \ref{exe.perm.cycles} \textbf{(a)} (applied to $1,2,\ldots,k$ instead
of $i_{1},i_{2},\ldots,i_{k}$) yields%
\[
\sigma\circ\operatorname*{cyc}\nolimits_{1,2,\ldots,k}\circ\sigma
^{-1}=\operatorname*{cyc}\nolimits_{\sigma\left(  1\right)  ,\sigma\left(
2\right)  ,\ldots,\sigma\left(  k\right)  }=\operatorname*{cyc}%
\nolimits_{i_{1},i_{2},\ldots,i_{k}}%
\]
(since $\left(  \sigma\left(  1\right)  ,\sigma\left(  2\right)
,\ldots,\sigma\left(  k\right)  \right)  =\left(  i_{1},i_{2},\ldots
,i_{k}\right)  $). Hence,
\[
\underbrace{\operatorname*{cyc}\nolimits_{i_{1},i_{2},\ldots,i_{k}}}%
_{=\sigma\circ\operatorname*{cyc}\nolimits_{1,2,\ldots,k}\circ\sigma^{-1}%
}\circ\sigma=\sigma\circ\operatorname*{cyc}\nolimits_{1,2,\ldots,k}%
\circ\underbrace{\sigma^{-1}\circ\sigma}_{=\operatorname*{id}}=\sigma
\circ\operatorname*{cyc}\nolimits_{1,2,\ldots,k}.
\]
Thus,%
\begin{align*}
\left(  -1\right)  ^{\operatorname*{cyc}\nolimits_{i_{1},i_{2},\ldots,i_{k}%
}\circ\sigma}  &  =\left(  -1\right)  ^{\sigma\circ\operatorname*{cyc}%
\nolimits_{1,2,\ldots,k}}=\left(  -1\right)  ^{\sigma}\cdot\underbrace{\left(
-1\right)  ^{\operatorname*{cyc}\nolimits_{1,2,\ldots,k}}}_{=\left(
-1\right)  ^{k-1}}\\
&  \ \ \ \ \ \ \ \ \ \ \left(  \text{by (\ref{eq.sign.prod}), applied to }%
\tau=\operatorname*{cyc}\nolimits_{1,2,\ldots,k}\right) \\
&  =\left(  -1\right)  ^{\sigma}\cdot\left(  -1\right)  ^{k-1}.
\end{align*}
Compared with%
\begin{align*}
\left(  -1\right)  ^{\operatorname*{cyc}\nolimits_{i_{1},i_{2},\ldots,i_{k}%
}\circ\sigma}  &  =\left(  -1\right)  ^{\operatorname*{cyc}\nolimits_{i_{1}%
,i_{2},\ldots,i_{k}}}\cdot\left(  -1\right)  ^{\sigma}\\
&  \ \ \ \ \ \ \ \ \ \ \left(  \text{by (\ref{eq.sign.prod}), applied to
}\operatorname*{cyc}\nolimits_{i_{1},i_{2},\ldots,i_{k}}\text{ and }%
\sigma\text{ instead of }\sigma\text{ and }\tau\right) \\
&  =\left(  -1\right)  ^{\sigma}\cdot\left(  -1\right)  ^{\operatorname*{cyc}%
\nolimits_{i_{1},i_{2},\ldots,i_{k}}},
\end{align*}
this yields $\left(  -1\right)  ^{\sigma}\cdot\left(  -1\right)
^{\operatorname*{cyc}\nolimits_{i_{1},i_{2},\ldots,i_{k}}}=\left(  -1\right)
^{\sigma}\cdot\left(  -1\right)  ^{k-1}$. We can cancel $\left(  -1\right)
^{\sigma}$ from this equality (since $\left(  -1\right)  ^{\sigma}\in\left\{
1,-1\right\}  $ is a nonzero integer), and thus obtain $\left(  -1\right)
^{\operatorname*{cyc}\nolimits_{i_{1},i_{2},\ldots,i_{k}}}=\left(  -1\right)
^{k-1}$. This solves Exercise \ref{exe.perm.cycles} \textbf{(d)}.
\end{vershort}

\begin{verlong}
Recall that $S_{n}$ is the set of all permutations of the set $\left\{
1,2,\ldots,n\right\}  $. In other words, $S_{n}$ is the set of all
permutations of the set $\left[  n\right]  $ (since $\left\{  1,2,\ldots
,n\right\}  =\left[  n\right]  $). Hence, $\sigma$ is a permutation of
$\left[  n\right]  $ (since $\sigma\in S_{n}$), and thus is a bijective map
from $\left[  n\right]  $ to $\left[  n\right]  $. Hence, the map $\sigma$ is
bijective and thus injective.

For every $p\in\left\{  1,2,\ldots,k\right\}  $, let $j_{p}$ be the element
$\sigma\left(  i_{p}\right)  \in\left[  n\right]  $. Then, $\left(
j_{1},j_{2},\ldots,j_{k}\right)  =\left(  \sigma\left(  i_{1}\right)
,\sigma\left(  i_{2}\right)  ,\ldots,\sigma\left(  i_{k}\right)  \right)  $.
Furthermore, $j_{1},j_{2},\ldots,j_{k}$ are $k$ distinct elements of $\left[
n\right]  $\ \ \ \ \footnote{\textit{Proof.} Clearly, $\sigma\left(
i_{1}\right)  ,\sigma\left(  i_{2}\right)  ,\ldots,\sigma\left(  i_{k}\right)
$ are $k$ elements of $\left[  n\right]  $ (since $\sigma$ is a map from
$\left[  n\right]  $ to $\left[  n\right]  $).
\par
The map $\sigma$ is injective. Hence, $\sigma$ sends distinct elements to
distinct elements. Thus, the elements $\sigma\left(  i_{1}\right)
,\sigma\left(  i_{2}\right)  ,\ldots,\sigma\left(  i_{k}\right)  $ are
distinct (since the elements $i_{1},i_{2},\ldots,i_{k}$ are distinct). Thus,
$\sigma\left(  i_{1}\right)  ,\sigma\left(  i_{2}\right)  ,\ldots
,\sigma\left(  i_{k}\right)  $ are $k$ distinct elements of $\left[  n\right]
$. In other words, $j_{1},j_{2},\ldots,j_{k}$ are $k$ distinct elements of
$\left[  n\right]  $ (since $\left(  j_{1},j_{2},\ldots,j_{k}\right)  =\left(
\sigma\left(  i_{1}\right)  ,\sigma\left(  i_{2}\right)  ,\ldots,\sigma\left(
i_{k}\right)  \right)  $).}. Therefore, $\operatorname*{cyc}\nolimits_{j_{1}%
,j_{2},\ldots,j_{k}}$ is a well-defined permutation in $S_{n}$.

We have defined $\operatorname*{cyc}\nolimits_{i_{1},i_{2},\ldots,i_{k}}$ to
be the permutation in $S_{n}$ which sends $i_{1},i_{2},\ldots,i_{k}$ to
$i_{2},i_{3},\ldots,i_{k},i_{1}$, respectively, while leaving all other
elements of $\left[  n\right]  $ fixed. Therefore:

\begin{itemize}
\item The permutation $\operatorname*{cyc}\nolimits_{i_{1},i_{2},\ldots,i_{k}%
}$ sends $i_{1},i_{2},\ldots,i_{k}$ to $i_{2},i_{3},\ldots,i_{k},i_{1}$,
respectively. In other words,%
\begin{equation}
\operatorname*{cyc}\nolimits_{i_{1},i_{2},\ldots,i_{k}}\left(  i_{p}\right)
=i_{p+1}\ \ \ \ \ \ \ \ \ \ \text{for every }p\in\left\{  1,2,\ldots
,k\right\}  , \label{sol.perm.cycles.a.cyc-manifest.1}%
\end{equation}
where $i_{k+1}$ means $i_{1}$.

\item The permutation $\operatorname*{cyc}\nolimits_{i_{1},i_{2},\ldots,i_{k}%
}$ leaves all other elements of $\left[  n\right]  $ fixed (where
\textquotedblleft other\textquotedblright\ means \textquotedblleft other than
$i_{1},i_{2},\ldots,i_{k}$\textquotedblright). In other words,%
\begin{equation}
\operatorname*{cyc}\nolimits_{i_{1},i_{2},\ldots,i_{k}}\left(  q\right)
=q\ \ \ \ \ \ \ \ \ \ \text{for every }q\in\left[  n\right]  \setminus\left\{
i_{1},i_{2},\ldots,i_{k}\right\}  . \label{sol.perm.cycles.a.cyc-manifest.2}%
\end{equation}

\end{itemize}

Furthermore:

\begin{itemize}
\item We have%
\begin{equation}
\operatorname*{cyc}\nolimits_{j_{1},j_{2},\ldots,j_{k}}\left(  j_{p}\right)
=j_{p+1}\ \ \ \ \ \ \ \ \ \ \text{for every }p\in\left\{  1,2,\ldots
,k\right\}  , \label{sol.perm.cycles.a.cyc-manifest.1'}%
\end{equation}
where $j_{k+1}$ means $j_{1}$. (This is proven in the same way as
(\ref{sol.perm.cycles.a.cyc-manifest.1}), except that every $i_{x}$ is
replaced by $j_{x}$.)

\item We have%
\begin{equation}
\operatorname*{cyc}\nolimits_{j_{1},j_{2},\ldots,j_{k}}\left(  q\right)
=q\ \ \ \ \ \ \ \ \ \ \text{for every }q\in\left[  n\right]  \setminus\left\{
j_{1},j_{2},\ldots,j_{k}\right\}  . \label{sol.perm.cycles.a.cyc-manifest.2'}%
\end{equation}
(This is proven in the same way as (\ref{sol.perm.cycles.a.cyc-manifest.2}),
except that every $i_{x}$ is replaced by $j_{x}$.)
\end{itemize}

In the following, we shall use the notation $i_{k+1}$ as a synonym for $i_{1}%
$, and the notation $j_{k+1}$ as a synonym for $j_{1}$. Then,%
\begin{equation}
j_{p}=\sigma\left(  i_{p}\right)  \ \ \ \ \ \ \ \ \ \ \text{for every }%
p\in\left\{  1,2,\ldots,k+1\right\}  \label{sol.perm.cycles.a.jp}%
\end{equation}
\footnote{\textit{Proof of (\ref{sol.perm.cycles.a.jp}):} Let $p\in\left\{
1,2,\ldots,k+1\right\}  $. We need to prove that $j_{p}=\sigma\left(
i_{p}\right)  $. If $p\in\left\{  1,2,\ldots,k\right\}  $, then this follows
from the definition of $j_{p}$. Hence, for the rest of this proof, we can WLOG
assume that we don't have $p\in\left\{  1,2,\ldots,k\right\}  $. Assume this.
\par
We have $p\in\left\{  1,2,\ldots,k+1\right\}  $, but we don't have
$p\in\left\{  1,2,\ldots,k\right\}  $. Hence, we have $p\in\left\{
1,2,\ldots,k+1\right\}  \setminus\left\{  1,2,\ldots,k\right\}  =\left\{
k+1\right\}  $. In other words, $p=k+1$, so that $i_{p}=i_{k+1}=i_{1}$ and
thus $\sigma\left(  i_{p}\right)  =\sigma\left(  i_{1}\right)  $. On the other
hand, from $p=k+1$, we obtain $j_{p}=j_{k+1}=j_{1}=\sigma\left(  i_{1}\right)
$ (by the definition of $j_{1}$). Compared with $\sigma\left(  i_{p}\right)
=\sigma\left(  i_{1}\right)  $, this yields $j_{p}=\sigma\left(  i_{p}\right)
$. Thus, $j_{p}=\sigma\left(  i_{p}\right)  $ is proven.}.

Now, let us show that
\begin{equation}
\left(  \sigma\circ\operatorname*{cyc}\nolimits_{i_{1},i_{2},\ldots,i_{k}%
}\circ\sigma^{-1}\right)  \left(  q\right)  =\operatorname*{cyc}%
\nolimits_{j_{1},j_{2},\ldots,j_{k}}\left(  q\right)
\label{sol.perm.cycles.a.qq}%
\end{equation}
for every $q\in\left[  n\right]  $.

[\textit{Proof of (\ref{sol.perm.cycles.a.qq}):} Let $q\in\left[  n\right]  $.
We must prove (\ref{sol.perm.cycles.a.qq}). We are in one of the following two cases:

\textit{Case 1:} We have $q\in\left\{  j_{1},j_{2},\ldots,j_{k}\right\}  $.

\textit{Case 2:} We have $q\notin\left\{  j_{1},j_{2},\ldots,j_{k}\right\}  $.

Let us first consider Case 1. In this case, we have $q\in\left\{  j_{1}%
,j_{2},\ldots,j_{k}\right\}  $. Thus, $q=j_{p}$ for some $p\in\left\{
1,2,\ldots,k\right\}  $. Consider this $p$. Clearly, $p+1\in\left\{
2,3,\ldots,k+1\right\}  \subseteq\left\{  1,2,\ldots,k+1\right\}  $. Hence,
applying (\ref{sol.perm.cycles.a.jp}) to $p+1$ instead of $p$, we obtain
$j_{p+1}=\sigma\left(  i_{p+1}\right)  $. But $q=j_{p}=\sigma\left(
i_{p}\right)  $ (by the definition of $j_{p}$) and thus $\sigma^{-1}\left(
q\right)  =i_{p}$. Hence,
\begin{align*}
\left(  \sigma\circ\operatorname*{cyc}\nolimits_{i_{1},i_{2},\ldots,i_{k}%
}\circ\sigma^{-1}\right)  \left(  q\right)   &  =\sigma\left(
\operatorname*{cyc}\nolimits_{i_{1},i_{2},\ldots,i_{k}}\left(
\underbrace{\sigma^{-1}\left(  q\right)  }_{=i_{p}}\right)  \right)
=\sigma\left(  \underbrace{\operatorname*{cyc}\nolimits_{i_{1},i_{2}%
,\ldots,i_{k}}\left(  i_{p}\right)  }_{\substack{=i_{p+1}\\\text{(by
(\ref{sol.perm.cycles.a.cyc-manifest.1}))}}}\right) \\
&  =\sigma\left(  i_{p+1}\right)  .
\end{align*}
Compared with%
\begin{align*}
\operatorname*{cyc}\nolimits_{j_{1},j_{2},\ldots,j_{k}}\left(  \underbrace{q}%
_{=j_{p}}\right)   &  =\operatorname*{cyc}\nolimits_{j_{1},j_{2},\ldots,j_{k}%
}\left(  j_{p}\right)  =j_{p+1}\ \ \ \ \ \ \ \ \ \ \left(  \text{by
(\ref{sol.perm.cycles.a.cyc-manifest.1'})}\right) \\
&  =\sigma\left(  i_{p+1}\right)  ,
\end{align*}
this yields $\left(  \sigma\circ\operatorname*{cyc}\nolimits_{i_{1}%
,i_{2},\ldots,i_{k}}\circ\sigma^{-1}\right)  \left(  q\right)
=\operatorname*{cyc}\nolimits_{j_{1},j_{2},\ldots,j_{k}}\left(  q\right)  $.
Thus, (\ref{sol.perm.cycles.a.qq}) is proven in Case 1.

Let us now consider Case 2. In this case, we have $q\notin\left\{  j_{1}%
,j_{2},\ldots,j_{k}\right\}  $. Hence, $\sigma^{-1}\left(  q\right)
\notin\left\{  i_{1},i_{2},\ldots,i_{k}\right\}  $%
\ \ \ \ \footnote{\textit{Proof.} Assume the contrary. Thus, $\sigma
^{-1}\left(  q\right)  \in\left\{  i_{1},i_{2},\ldots,i_{k}\right\}  $. In
other words, there exists a $p\in\left\{  1,2,\ldots,k\right\}  $ such that
$\sigma^{-1}\left(  q\right)  =i_{p}$. Consider this $p$. We have $\sigma
^{-1}\left(  q\right)  =i_{p}$, thus $q=\sigma\left(  i_{p}\right)  =j_{p}$
(since $j_{p}$ is defined as $\sigma\left(  i_{p}\right)  $). Thus,
$q=j_{p}\in\left\{  j_{1},j_{2},\ldots,j_{k}\right\}  $; but this contradicts
$q\notin\left\{  j_{1},j_{2},\ldots,j_{k}\right\}  $. This contradiction shows
that our assumption was wrong, qed.}. Combining this with $\sigma^{-1}\left(
q\right)  \in\left[  n\right]  $, we obtain $\sigma^{-1}\left(  q\right)
\in\left[  n\right]  \setminus\left\{  i_{1},i_{2},\ldots,i_{k}\right\}  $.
Therefore, (\ref{sol.perm.cycles.a.cyc-manifest.2}) (applied to $\sigma
^{-1}\left(  q\right)  $ instead of $q$) yields $\operatorname*{cyc}%
\nolimits_{i_{1},i_{2},\ldots,i_{k}}\left(  \sigma^{-1}\left(  q\right)
\right)  =\sigma^{-1}\left(  q\right)  $. Hence,%
\begin{equation}
\left(  \sigma\circ\operatorname*{cyc}\nolimits_{i_{1},i_{2},\ldots,i_{k}%
}\circ\sigma^{-1}\right)  \left(  q\right)  =\sigma\left(
\underbrace{\operatorname*{cyc}\nolimits_{i_{1},i_{2},\ldots,i_{k}}\left(
\sigma^{-1}\left(  q\right)  \right)  }_{=\sigma^{-1}\left(  q\right)
}\right)  =\sigma\left(  \sigma^{-1}\left(  q\right)  \right)  =q.
\label{sol.perm.cycles.a.qq.pf.5}%
\end{equation}
On the other hand, $q\in\left[  n\right]  \setminus\left\{  j_{1},j_{2}%
,\ldots,j_{k}\right\}  $ (since $q\in\left[  n\right]  $ and $q\notin\left\{
j_{1},j_{2},\ldots,j_{k}\right\}  $) and thus $\operatorname*{cyc}%
\nolimits_{j_{1},j_{2},\ldots,j_{k}}\left(  q\right)  =q$ (by
(\ref{sol.perm.cycles.a.cyc-manifest.2'})). Compared with
(\ref{sol.perm.cycles.a.qq.pf.5}), this yields \newline$\left(  \sigma
\circ\operatorname*{cyc}\nolimits_{i_{1},i_{2},\ldots,i_{k}}\circ\sigma
^{-1}\right)  \left(  q\right)  =\operatorname*{cyc}\nolimits_{j_{1}%
,j_{2},\ldots,j_{k}}\left(  q\right)  $. Thus, (\ref{sol.perm.cycles.a.qq}) is
proven in Case 2.

We have now proven (\ref{sol.perm.cycles.a.qq}) in each of the two Cases 1 and
2. Therefore, (\ref{sol.perm.cycles.a.qq}) always holds.]

Thus, we know that $\left(  \sigma\circ\operatorname*{cyc}\nolimits_{i_{1}%
,i_{2},\ldots,i_{k}}\circ\sigma^{-1}\right)  \left(  q\right)
=\operatorname*{cyc}\nolimits_{j_{1},j_{2},\ldots,j_{k}}\left(  q\right)  $
for every $q\in\left[  n\right]  $. In other words, $\sigma\circ
\operatorname*{cyc}\nolimits_{i_{1},i_{2},\ldots,i_{k}}\circ\sigma
^{-1}=\operatorname*{cyc}\nolimits_{j_{1},j_{2},\ldots,j_{k}}$ (because both
$\sigma\circ\operatorname*{cyc}\nolimits_{i_{1},i_{2},\ldots,i_{k}}\circ
\sigma^{-1}$ and $\operatorname*{cyc}\nolimits_{j_{1},j_{2},\ldots,j_{k}}$ are
maps from $\left[  n\right]  $ to $\left[  n\right]  $). Hence,%
\[
\sigma\circ\operatorname*{cyc}\nolimits_{i_{1},i_{2},\ldots,i_{k}}\circ
\sigma^{-1}=\operatorname*{cyc}\nolimits_{j_{1},j_{2},\ldots,j_{k}%
}=\operatorname*{cyc}\nolimits_{\sigma\left(  i_{1}\right)  ,\sigma\left(
i_{2}\right)  ,\ldots,\sigma\left(  i_{k}\right)  }%
\]
(since $\left(  j_{1},j_{2},\ldots,j_{k}\right)  =\left(  \sigma\left(
i_{1}\right)  ,\sigma\left(  i_{2}\right)  ,\ldots,\sigma\left(  i_{k}\right)
\right)  $). This solves Exercise \ref{exe.perm.cycles} \textbf{(a)}.

\textbf{(b)} Let $p\in\left\{  0,1,\ldots,n-k\right\}  $. Then, $p\geq0$ and
$p\leq n-k$. Hence,
\[
\left\{  p+1,p+2,\ldots,p+k\right\}  \subseteq\left\{  1,2,\ldots,n\right\}
\]
(since $\underbrace{p}_{\geq0}+1\geq1$ and $\underbrace{p}_{\leq n-k}+k\leq
n-k+k=n$). Thus, the elements $p+1,p+2,\ldots,p+k$ belong to $\left\{
1,2,\ldots,n\right\}  =\left[  n\right]  $. Hence, $\operatorname*{cyc}%
\nolimits_{p+1,p+2,\ldots,p+k}$ is well-defined. We let $\sigma$ denote the
permutation $\operatorname*{cyc}\nolimits_{p+1,p+2,\ldots,p+k}$.

The permutation $\operatorname*{cyc}\nolimits_{p+1,p+2,\ldots,p+k}$ is defined
to be the permutation in $S_{n}$ which sends $p+1,p+2,\ldots,p+k$ to
$p+2,p+3,\ldots,p+k,p+1$, respectively, while leaving all other elements of
$\left[  n\right]  $ fixed. Therefore:

\begin{itemize}
\item The permutation $\operatorname*{cyc}\nolimits_{p+1,p+2,\ldots,p+k}$
sends $p+1,p+2,\ldots,p+k$ to $p+2,p+3,\ldots,p+k,p+1$, respectively. In other
words, the permutation $\sigma$ sends $p+1,p+2,\ldots,p+k$ to $p+2,p+3,\ldots
,p+k,p+1$, respectively (since $\sigma=\operatorname*{cyc}%
\nolimits_{p+1,p+2,\ldots,p+k}$). In other words, we have%
\begin{equation}
\left(  \sigma\left(  p+i\right)  =p+\left(  i+1\right)
\ \ \ \ \ \ \ \ \ \ \text{for every }i\in\left\{  1,2,\ldots,k-1\right\}
\right)  \label{sol.perm.cycles.b.sigma-manifest.1a}%
\end{equation}
and%
\begin{equation}
\sigma\left(  p+k\right)  =p+1. \label{sol.perm.cycles.b.sigma-manifest.1b}%
\end{equation}

\item The permutation $\operatorname*{cyc}\nolimits_{p+1,p+2,\ldots,p+k}$
leaves all other elements of $\left[  n\right]  $ fixed (where
\textquotedblleft other\textquotedblright\ means \textquotedblleft other than
$p+1,p+2,\ldots,p+k$\textquotedblright). In other words,
\[
\operatorname*{cyc}\nolimits_{p+1,p+2,\ldots,p+k}\left(  q\right)
=q\ \ \ \ \ \ \ \ \ \ \text{for every }q\in\left[  n\right]  \setminus\left\{
p+1,p+2,\ldots,p+k\right\}  .
\]
In other words,%
\begin{equation}
\sigma\left(  q\right)  =q\ \ \ \ \ \ \ \ \ \ \text{for every }q\in\left[
n\right]  \setminus\left\{  p+1,p+2,\ldots,p+k\right\}
\label{sol.perm.cycles.b.sigma-manifest.2}%
\end{equation}
(since $\sigma=\operatorname*{cyc}\nolimits_{p+1,p+2,\ldots,p+k}$).
\end{itemize}

We can now observe that
\begin{equation}
\sigma\left(  q\right)  \geq q\ \ \ \ \ \ \ \ \ \ \text{for every }q\in\left[
n\right]  \text{ satisfying }q\neq p+k \label{sol.perm.cycles.b.sigma.1}%
\end{equation}
\footnote{\textit{Proof of (\ref{sol.perm.cycles.b.sigma.1}):} Let
$q\in\left[  n\right]  $ be such that $q\neq p+k$. We must prove that
$\sigma\left(  q\right)  \geq q$.
\par
We are in one of the following two cases:
\par
\textit{Case 1:} We have $q\in\left\{  p+1,p+2,\ldots,p+k\right\}  $.
\par
\textit{Case 2:} We have $q\notin\left\{  p+1,p+2,\ldots,p+k\right\}  $.
\par
Let us first consider Case 1. In this case, we have $q\in\left\{
p+1,p+2,\ldots,p+k\right\}  $. Hence, $q=p+i$ for some $i\in\left\{
1,2,\ldots,k\right\}  $. Consider this $i$. We have $q\neq p+k$, so that
$q-p\neq k$ and thus $k\neq q-p=i$ (since $q=p+i$). Therefore, $i\neq k$.
Combined with $i\in\left\{  1,2,\ldots,k\right\}  $, this shows that
$i\in\left\{  1,2,\ldots,k\right\}  \setminus\left\{  k\right\}  =\left\{
1,2,\ldots,k-1\right\}  $. Hence, (\ref{sol.perm.cycles.b.sigma-manifest.1a})
yields $\sigma\left(  p+i\right)  =p+\underbrace{\left(  i+1\right)  }_{\geq
i}\geq p+i=q$. Hence, $\sigma\left(  \underbrace{q}_{=p+i}\right)
=\sigma\left(  p+i\right)  \geq q$. We thus have proven $\sigma\left(
q\right)  \geq q$ in Case 1.
\par
Let us now consider Case 2. In this case, we have $q\notin\left\{
p+1,p+2,\ldots,p+k\right\}  $. Hence, $q\in\left[  n\right]  \setminus\left\{
p+1,p+2,\ldots,p+k\right\}  $. Therefore, $\sigma\left(  q\right)  =q$ (by
(\ref{sol.perm.cycles.b.sigma-manifest.2})). Thus, $\sigma\left(  q\right)
\geq q$ is proven in Case 2.
\par
We have now proven $\sigma\left(  q\right)  \geq q$ in both Cases 1 and 2.
Therefore, $\sigma\left(  q\right)  \geq q$ always holds. This proves
(\ref{sol.perm.cycles.b.sigma.1}).}. Furthermore,%
\begin{equation}
\sigma\left(  q\right)  \leq q+1\ \ \ \ \ \ \ \ \ \ \text{for every }%
q\in\left[  n\right]  \label{sol.perm.cycles.b.sigma.2}%
\end{equation}
\footnote{\textit{Proof of (\ref{sol.perm.cycles.b.sigma.2}):} Let
$q\in\left[  n\right]  $. We must prove that $\sigma\left(  q\right)  \leq
q+1$.
\par
Assume the contrary (for the sake of contradiction). Thus, $\sigma\left(
q\right)  >q+1$. If we had $q\in\left[  n\right]  \setminus\left\{
p+1,p+2,\ldots,p+k\right\}  $, then we would have $\sigma\left(  q\right)  =q$
(by (\ref{sol.perm.cycles.b.sigma-manifest.2})), which would contradict
$\sigma\left(  q\right)  >q+1>q$. Thus, we cannot have $q\in\left[  n\right]
\setminus\left\{  p+1,p+2,\ldots,p+k\right\}  $. We therefore have%
\[
q\in\left[  n\right]  \setminus\left(  \left[  n\right]  \setminus\left\{
p+1,p+2,\ldots,p+k\right\}  \right)  \subseteq\left\{  p+1,p+2,\ldots
,p+k\right\}  .
\]
Hence, $q=p+i$ for some $i\in\left\{  1,2,\ldots,k\right\}  $. Consider this
$i$. Clearly, $i\geq1$ and $i\leq k$.
\par
If we had $i=k$, then we would have%
\begin{align*}
\sigma\left(  \underbrace{q}_{=p+i}\right)   &  =\sigma\left(
p+\underbrace{i}_{=k}\right)  =\sigma\left(  p+k\right)  =p+\underbrace{1}%
_{\leq i}\ \ \ \ \ \ \ \ \ \ \left(  \text{by
(\ref{sol.perm.cycles.b.sigma-manifest.1b})}\right) \\
&  \leq p+i=q,
\end{align*}
which would contradict $\sigma\left(  q\right)  >q+1>q$. Thus, we cannot have
$i=k$. We thus have $i\neq k$. Combined with $i\in\left\{  1,2,\ldots
,k\right\}  $, this shows that $i\in\left\{  1,2,\ldots,k\right\}
\setminus\left\{  k\right\}  =\left\{  1,2,\ldots,k-1\right\}  $. Hence,
(\ref{sol.perm.cycles.b.sigma-manifest.1a}) yields $\sigma\left(  p+i\right)
=p+\left(  i+1\right)  =\underbrace{p+i}_{=q}+1=q+1$. This contradicts
$\sigma\left(  \underbrace{p+i}_{=q}\right)  =\sigma\left(  q\right)  >q+1$.
This contradiction proves that our assumption was wrong. Hence, $\sigma\left(
q\right)  \leq q+1$ is proven. We are thus done proving
(\ref{sol.perm.cycles.b.sigma.2}).}.

Now, set%
\[
A=\left\{  \left(  p+h,p+k\right)  \ \mid\ h\in\left\{  1,2,\ldots
,k-1\right\}  \right\}  .
\]
In other words,%
\[
A=\left\{  \left(  p+1,p+k\right)  ,\left(  p+2,p+k\right)  ,\ldots,\left(
p+\left(  k-1\right)  ,p+k\right)  \right\}  .
\]
Thus, the set $A$ has $k-1$ elements (since the $k-1$ pairs \newline$\left(
p+1,p+k\right)  ,\left(  p+2,p+k\right)  ,\ldots,\left(  p+\left(  k-1\right)
,p+k\right)  $ are clearly distinct). In other words, $\left\vert A\right\vert
=k-1$.

Now, let $\operatorname*{Inv}\sigma$ denote the set of all inversions of
$\sigma$. Then, $\ell\left(  \sigma\right)  =\left\vert \operatorname*{Inv}%
\sigma\right\vert $\ \ \ \ \footnote{\textit{Proof.} Recall that $\ell\left(
\sigma\right)  $ is defined as the number of inversions of $\sigma$. Thus,%
\[
\ell\left(  \sigma\right)  =\left(  \text{the number of inversions of }%
\sigma\right)  =\left\vert \underbrace{\left(  \text{the set of all inversions
of }\sigma\right)  }_{=\operatorname*{Inv}\sigma}\right\vert =\left\vert
\operatorname*{Inv}\sigma\right\vert ,
\]
qed.}.

But $A\subseteq\operatorname*{Inv}\sigma$\ \ \ \ \footnote{\textit{Proof.} Let
$c\in A$. We shall show that $c\in\operatorname*{Inv}\sigma$.
\par
We have $c\in A=\left\{  \left(  p+h,p+k\right)  \ \mid\ h\in\left\{
1,2,\ldots,k-1\right\}  \right\}  $. Hence, $c$ can be written in the form
$c=\left(  p+h,p+k\right)  $ for some $h\in\left\{  1,2,\ldots,k-1\right\}  $.
Consider this $h$. From $h\in\left\{  1,2,\ldots,k-1\right\}  $, we obtain
$1\leq h\leq k-1$, so that $h\leq k-1<k$. Thus, $p+\underbrace{h}_{<k}<p+k$.
Moreover, $1\leq p+h$ (since $\underbrace{p}_{\geq0}+\underbrace{h}_{\geq
1}\geq0+1=1$) and $p+k\leq n$ (since $p\leq n-k$). Thus, $1\leq p+h<p+k\leq
n$.
\par
Applying (\ref{sol.perm.cycles.b.sigma-manifest.1a}) to $i=h$, we obtain
$\sigma\left(  p+h\right)  =p+\left(  \underbrace{h}_{\geq1>0}+1\right)
>p+\left(  0+1\right)  =p+1=\sigma\left(  p+k\right)  $ (by
(\ref{sol.perm.cycles.b.sigma-manifest.1b})).
\par
Now, $\left(  p+h,p+k\right)  $ is a pair of integers satisfying $1\leq
p+h<p+k\leq n$ and $\sigma\left(  p+h\right)  >\sigma\left(  p+k\right)  $. In
other words, $\left(  p+h,p+k\right)  $ is a pair of integers $\left(
i,j\right)  $ satisfying $1\leq i<j\leq n$ and $\sigma\left(  i\right)
>\sigma\left(  j\right)  $. In other words, $\left(  p+h,p+k\right)  $ is an
inversion of $\sigma$ (by the definition of \textquotedblleft inversion of
$\sigma$\textquotedblright). In other words, $\left(  p+h,p+k\right)
\in\operatorname*{Inv}\sigma$ (since $\operatorname*{Inv}\sigma$ is the set of
all inversions of $\sigma$). Thus, $c=\left(  p+h,p+k\right)  \in
\operatorname*{Inv}\sigma$.
\par
Now, let us forget that we fixed $c$. We thus have shown that $c\in
\operatorname*{Inv}\sigma$ for every $c\in A$. In other words, $A\subseteq
\operatorname*{Inv}\sigma$, qed.} and $\operatorname*{Inv}\sigma\subseteq
A$\ \ \ \ \footnote{\textit{Proof.} Let $c\in\operatorname*{Inv}\sigma$. We
shall show that $c\in A$.
\par
We have $c\in\operatorname*{Inv}\sigma$. In other words, $c$ is an inversion
of $\sigma$ (since $\operatorname*{Inv}\sigma$ is the set of all inversions of
$\sigma$). In other words, $c$ is a pair $\left(  i,j\right)  $ of integers
satisfying $1\leq i<j\leq n$ and $\sigma\left(  i\right)  >\sigma\left(
j\right)  $. Consider this $\left(  i,j\right)  $. Thus, $c=\left(
i,j\right)  $.
\par
We have $\sigma\left(  i\right)  >\sigma\left(  j\right)  $, so that
$\sigma\left(  i\right)  \geq\sigma\left(  j\right)  +1$ (since $\sigma\left(
i\right)  $ and $\sigma\left(  j\right)  $ are integers). Also, $i<j$, so that
$i\leq j-1$ (since $i$ and $j$ are integers). In other words, $i+1\leq j$. But
(\ref{sol.perm.cycles.b.sigma.2}) (applied to $q=i$) yields $\sigma\left(
i\right)  \leq i+1\leq j$.
\par
Let us first show that $j=p+k$. Indeed, let us assume the contrary (for the
sake of contradiction). Thus, $j\neq p+k$. Hence, $\sigma\left(  j\right)
\geq j$ (by (\ref{sol.perm.cycles.b.sigma.1}), applied to $q=j$). Now,
$\sigma\left(  i\right)  >\sigma\left(  j\right)  \geq j$. This contradicts
$\sigma\left(  i\right)  \leq j$. This contradiction shows that our assumption
was wrong. Hence, $j=p+k$ is proven.
\par
Now, $\sigma\left(  \underbrace{j}_{=p+k}\right)  =\sigma\left(  p+k\right)
=p+1$ (by (\ref{sol.perm.cycles.b.sigma-manifest.1b})). Hence, $\sigma\left(
i\right)  >\sigma\left(  j\right)  =p+1$. Therefore, $p+1<\sigma\left(
i\right)  \leq i+1$. Subtracting $1$ from both sides of this inequality, we
obtain $p<i$. Hence, $i>p$, so that $i\geq p+1$ (since $i$ and $p$ are
integers). Combined with $i<j=p+k$, this yields $i\in\left\{  p+1,p+2,\ldots
,p+k-1\right\}  $. Thus, $i-p\in\left\{  1,2,\ldots,k-1\right\}  $.
\par
So we know that the element $i-p\in\left\{  1,2,\ldots,k-1\right\}  $
satisfies $c=\left(  \underbrace{i}_{=p+\left(  i-p\right)  },\underbrace{j}%
_{=p+k}\right)  =\left(  p+\left(  i-p\right)  ,p+k\right)  $. Hence, there
exists an $h\in\left\{  1,2,\ldots,k-1\right\}  $ such that $c=\left(
p+h,p+k\right)  $ (namely, $h=i-p$). Thus,%
\[
c\in\left\{  \left(  p+h,p+k\right)  \ \mid\ h\in\left\{  1,2,\ldots
,k-1\right\}  \right\}  =A.
\]
\par
Now, let us forget that we fixed $c$. We thus have proven that $c\in A$ for
every $c\in\operatorname*{Inv}\sigma$. In other words, $\operatorname*{Inv}%
\sigma\subseteq A$, qed.}. Combining these two relations, we obtain
$A=\operatorname*{Inv}\sigma$. Hence, $\left\vert A\right\vert =\left\vert
\operatorname*{Inv}\sigma\right\vert $. Compared with $\ell\left(
\sigma\right)  =\left\vert \operatorname*{Inv}\sigma\right\vert $, this yields
$\ell\left(  \sigma\right)  =\left\vert A\right\vert =k-1$. This rewrites as
$\ell\left(  \operatorname*{cyc}\nolimits_{p+1,p+2,\ldots,p+k}\right)  =k-1$
(since $\sigma=\operatorname*{cyc}\nolimits_{p+1,p+2,\ldots,p+k}$). This
solves Exercise \ref{exe.perm.cycles} \textbf{(b)}.

\textbf{(c)} We shall only sketch the proof, since we will not use the result
in the following.

Let $i_{1},i_{2},\ldots,i_{k}$ be $k$ distinct elements of $\left[  n\right]
$. We extend the $k$-tuple $\left(  i_{1},i_{2},\ldots,i_{k}\right)  $ to an
infinite sequence $\left(  i_{1},i_{2},i_{3},\ldots\right)  $ of elements of
$\left[  n\right]  $ by setting%
\[
\left(  i_{u}=i_{\left(  \text{the element }u^{\prime}\in\left\{
1,2,\ldots,k\right\}  \text{ satisfying }u^{\prime}\equiv u\operatorname{mod}%
k\right)  }\ \ \ \ \ \ \ \ \ \ \text{for every }u\geq1\right)  .
\]
This sequence $\left(  i_{1},i_{2},i_{3},\ldots\right)  $ is periodic with
period $k$. In other words,%
\begin{equation}
i_{u}=i_{u+k}\ \ \ \ \ \ \ \ \ \ \text{for every }u\geq1.
\label{sol.perm.cycles.c.period}%
\end{equation}
From this, it is easy to obtain that%
\[
i_{1}+i_{2}+\cdots+i_{k}=i_{2}+i_{3}+\cdots+i_{k+1}=i_{3}+i_{4}+\cdots
+i_{k+2}=\cdots
\]
\footnote{\textit{Proof.} Every $u\in\left\{  1,2,3,\ldots\right\}  $
satisfies%
\begin{align*}
i_{u}+i_{u+1}+\cdots+i_{u+k-1}  &  =\underbrace{i_{u}}_{\substack{=i_{u+k}%
\\\text{(by (\ref{sol.perm.cycles.c.period}))}}}+\left(  i_{u+1}%
+i_{u+2}+\cdots+i_{u+k-1}\right) \\
&  =i_{u+k}+\left(  i_{u+1}+i_{u+2}+\cdots+i_{u+k-1}\right)  =\left(
i_{u+1}+i_{u+2}+\cdots+i_{u+k-1}\right)  +i_{u+k}\\
&  =i_{u+1}+i_{u+2}+\cdots+i_{u+k}.
\end{align*}
Thus, $i_{1}+i_{2}+\cdots+i_{k}=i_{2}+i_{3}+\cdots+i_{k+1}=i_{3}+i_{4}%
+\cdots+i_{k+2}=\cdots$, qed.}. Thus,%
\begin{equation}
i_{r+1}+i_{r+2}+\cdots+i_{r+k}=i_{1}+i_{2}+\cdots+i_{k}%
\ \ \ \ \ \ \ \ \ \ \text{for every }r\in\mathbb{N}.
\label{sol.perm.cycles.c.periodsum}%
\end{equation}

Let $\sigma=\operatorname*{cyc}\nolimits_{i_{1},i_{2},\ldots,i_{k}}$. Let
$\operatorname*{Inv}\sigma$ denote the set of all inversions of $\sigma$.
Then, $\ell\left(  \sigma\right)  =\left\vert \operatorname*{Inv}%
\sigma\right\vert $. (This can be seen as in the solution to Exercise
\ref{exe.perm.cycles} \textbf{(b)}.) Moreover, the definitions of the sequence
$\left(  i_{1},i_{2},i_{3},\ldots\right)  $ and of $\sigma$ show that%
\begin{equation}
\sigma\left(  i_{p}\right)  =i_{p+1}\ \ \ \ \ \ \ \ \ \ \text{for every }%
p\geq1. \label{sol.perm.cycles.c.sigma}%
\end{equation}

Now, fix $r\in\left\{  1,2,\ldots,k-1\right\}  $. We shall prove that
\begin{equation}
\text{there exists some }u\geq1\text{ such that }\left(  i_{u},i_{u+r}\right)
\in\operatorname*{Inv}\sigma. \label{sol.perm.cycles.c.mainclaim}%
\end{equation}

[\textit{Proof of (\ref{sol.perm.cycles.c.mainclaim}):} The sequence $\left(
i_{1},i_{2},i_{3},\ldots\right)  $ is periodic with period $k$, and its first
$k$ entries $i_{1},i_{2},\ldots,i_{k}$ are distinct. Hence, each entry of this
sequence repeats itself each $k$ steps, but not more often. Hence, every
integer $u\geq1$ satisfies $i_{u+r}\neq i_{u}$ (since $r\in\left\{
1,2,\ldots,k-1\right\}  $). In other words, every $u\geq1$ satisfies%
\begin{equation}
i_{u+r}-i_{u}\neq0. \label{sol.perm.cycles.c.mainclaim.pf.4}%
\end{equation}
The $k$-tuple $\left(  i_{1+r}-i_{1},i_{2+r}-i_{2},\ldots,i_{k+r}%
-i_{k}\right)  $ contains at least one positive entry\footnote{\textit{Proof.}
Assume the contrary. Thus, the $k$-tuple $\left(  i_{1+r}-i_{1},i_{2+r}%
-i_{2},\ldots,i_{k+r}-i_{k}\right)  $ contains no positive entries. In other
words, no $u\in\left\{  1,2,\ldots,k\right\}  $ satisfies $i_{u+r}-i_{u}>0$.
In other words, every $u\in\left\{  1,2,\ldots,k\right\}  $ satisfies
$i_{u+r}-i_{u}\leq0$. But (\ref{sol.perm.cycles.c.mainclaim.pf.4}) shows that
every $u\in\left\{  1,2,\ldots,k\right\}  $ satisfies $i_{u+r}-i_{u}\neq0$.
Thus, every $u\in\left\{  1,2,\ldots,k\right\}  $ satisfies $i_{u+r}-i_{u}<0$
(since $i_{u+r}-i_{u}\leq0$ and $i_{u+r}-i_{u}\neq0$). Hence,%
\[
\sum_{u=1}^{k}\underbrace{\left(  i_{u+r}-i_{u}\right)  }_{<0}<\sum_{u=1}%
^{k}0=0.
\]
But this contradicts%
\begin{align*}
\sum_{u=1}^{k}\left(  \underbrace{i_{u+r}}_{=i_{r+u}}-i_{u}\right)   &
=\sum_{u=1}^{k}\left(  i_{r+u}-i_{u}\right)  =\left(  \sum_{u=1}^{k}%
i_{r+u}\right)  -\left(  \sum_{u=1}^{k}i_{u}\right) \\
&  =\left(  i_{r+1}+i_{r+2}+\cdots+i_{r+k}\right)  -\left(  i_{1}+i_{2}%
+\cdots+i_{k}\right)  =0\ \ \ \ \ \ \ \ \ \ \left(  \text{by
(\ref{sol.perm.cycles.c.periodsum})}\right)  .
\end{align*}
This contradiction shows that our assumption was wrong, qed.}, and at least
one negative entry\footnote{This is proven similarly.}. Hence, there exists at
least one $u\geq1$ such that $i_{u+r}-i_{u}>0$ but $i_{\left(  u+1\right)
+r}-i_{u+1}<0$\ \ \ \ \footnote{\textit{Proof.} Assume the contrary. Then,
there exists no $u\geq1$ such that $i_{u+r}-i_{u}>0$ but $i_{\left(
u+1\right)  +r}-i_{u+1}<0$. Hence,%
\begin{equation}
\text{every }u\geq1\text{ satisfying }i_{u+r}-i_{u}>0\text{ must satisfy
}i_{\left(  u+1\right)  +r}-i_{u+1}\geq0.
\label{sol.perm.cycles.c.mainclaim.pf.5o}%
\end{equation}
Therefore,%
\begin{equation}
\text{every }u\geq1\text{ satisfying }i_{u+r}-i_{u}>0\text{ must satisfy
}i_{\left(  u+1\right)  +r}-i_{u+1}>0 \label{sol.perm.cycles.c.mainclaim.pf.5}%
\end{equation}
(because (\ref{sol.perm.cycles.c.mainclaim.pf.5o}) shows that $i_{\left(
u+1\right)  +r}-i_{u+1}\geq0$; but combining this with $i_{\left(  u+1\right)
+r}-i_{u+1}\neq0$ (which follows from (\ref{sol.perm.cycles.c.mainclaim.pf.4}%
), applied to $u+1$ instead of $u$), we obtain $i_{\left(  u+1\right)
+r}-i_{u+1}>0$).
\par
But there exists some $v\in\left\{  1,2,\ldots,k\right\}  $ satisfying
$i_{v+r}-i_{v}>0$ (since the $k$-tuple $\left(  i_{1+r}-i_{1},i_{2+r}%
-i_{2},\ldots,i_{k+r}-i_{k}\right)  $ contains at least one positive entry).
Consider this $v$. Then, we have $i_{v+r}-i_{v}>0$, therefore $i_{\left(
v+1\right)  +r}-i_{v+1}>0$ (by (\ref{sol.perm.cycles.c.mainclaim.pf.5}),
applied to $u=v$), therefore $i_{\left(  v+2\right)  +r}-i_{v+2}>0$ (by
(\ref{sol.perm.cycles.c.mainclaim.pf.5}), applied to $u=v+1$), therefore
$i_{\left(  v+3\right)  +r}-i_{v+3}>0$ (by
(\ref{sol.perm.cycles.c.mainclaim.pf.5}), applied to $u=v+2$), and so on.
Altogether, we thus obtain
\[
i_{h+r}-i_{h}>0\ \ \ \ \ \ \ \ \ \ \text{for every }h\geq v.
\]
In other words, $i_{h}<i_{h+r}$ for every $h\geq v$. Hence, $i_{v}%
<i_{v+r}<i_{v+2r}<i_{v+3r}<\cdots$. Thus, the numbers $i_{v},i_{v+r}%
,i_{v+2r},i_{v+3r},\ldots$ are pairwise distinct; hence, the sequence $\left(
i_{1},i_{2},i_{3},\ldots\right)  $ contains infinitely many distinct entries.
But this contradicts the fact that this sequence is periodic. This
contradiction proves that our assumption was wrong, qed.}. Consider this $u$.
We have $i_{u}<i_{u+r}$ (since $i_{u+r}-i_{u}>0$), so that $1\leq
i_{u}<i_{u+r}\leq n$. Also, (\ref{sol.perm.cycles.c.sigma}) (applied to $p=u$)
yields $\sigma\left(  i_{u}\right)  =i_{u+1}$. Moreover,
(\ref{sol.perm.cycles.c.sigma}) (applied to $p=u+r$) yields $\sigma\left(
i_{u+r}\right)  =i_{u+r+1}=i_{\left(  u+1\right)  +r}$. Hence,
\begin{align*}
\sigma\left(  i_{u}\right)   &  =i_{u+1}>i_{\left(  u+1\right)  +r}%
\ \ \ \ \ \ \ \ \ \ \left(  \text{since }i_{\left(  u+1\right)  +r}%
-i_{u+1}<0\right) \\
&  =\sigma\left(  i_{u+r}\right)  .
\end{align*}
So we know that $\left(  i_{u},i_{u+r}\right)  $ is a pair of integers
satisfying $1\leq i_{u}<i_{u+r}\leq n$ and $\sigma\left(  i_{u}\right)
>\sigma\left(  i_{u+r}\right)  $. In other words, $\left(  i_{u}%
,i_{u+r}\right)  $ is an inversion of $\sigma$ (by the definition of an
\textquotedblleft inversion\textquotedblright). In other words, $\left(
i_{u},i_{u+r}\right)  \in\operatorname*{Inv}\sigma$. Thus, we have found a
$u\geq1$ such that $\left(  i_{u},i_{u+r}\right)  \in\operatorname*{Inv}%
\sigma$. This proves (\ref{sol.perm.cycles.c.mainclaim}).]

Now, let us forget that we fixed $r$. We have shown that, for every
$r\in\left\{  1,2,\ldots,k-1\right\}  $, there exists some $u\geq1$ such that
$\left(  i_{u},i_{u+r}\right)  \in\operatorname*{Inv}\sigma$. Let us denote
this $u$ by $u_{r}$. Therefore, for every $r\in\left\{  1,2,\ldots
,k-1\right\}  $, we have found a $u_{r}\geq1$ such that $\left(  i_{u_{r}%
},i_{u_{r}+r}\right)  \in\operatorname*{Inv}\sigma$. The $k-1$ pairs%
\[
\left(  i_{u_{1}},i_{u_{1}+1}\right)  ,\ \left(  i_{u_{2}},i_{u_{2}+2}\right)
,\ \ldots,\ \left(  i_{u_{k-1}},i_{u_{k-1}+\left(  k-1\right)  }\right)
\]
are pairwise distinct\footnote{\textit{Proof.} Assume the contrary. Then,
there exist two distinct elements $x$ and $y$ of $\left\{  1,2,\ldots
,k-1\right\}  $ such that $\left(  i_{u_{x}},i_{u_{x}+x}\right)  =\left(
i_{u_{y}},i_{u_{y}+y}\right)  $. Consider these $x$ and $y$.
\par
We have $\left(  i_{u_{x}},i_{u_{x}+x}\right)  =\left(  i_{u_{y}},i_{u_{y}%
+y}\right)  $. In other words, $i_{u_{x}}=i_{u_{y}}$ and $i_{u_{x}+x}%
=i_{u_{y}+y}$. Since the numbers $i_{1},i_{2},\ldots,i_{k}$ are distinct (and
the sequence $\left(  i_{1},i_{2},i_{3},\ldots\right)  $ consists of these
numbers, repeated over and over), we obtain $u_{x}\equiv u_{y}%
\operatorname{mod}k$ from $i_{u_{x}}=i_{u_{y}}$, and we obtain $u_{x}+x\equiv
u_{y}+y\operatorname{mod}k$ from $i_{u_{x}+x}=i_{u_{y}+y}$. Subtracting the
congruence $u_{x}\equiv u_{y}\operatorname{mod}k$ from the congruence
$u_{x}+x\equiv u_{y}+y\operatorname{mod}k$, we obtain $x\equiv
y\operatorname{mod}k$. In light of $x,y\in\left\{  1,2,\ldots,k-1\right\}  $,
this shows that $x=y$. But this contradicts the fact that $x$ and $y$ are
distinct. This contradiction proves that our assumption was wrong, qed.}, and
all of them belong to $\operatorname*{Inv}\sigma$. Hence, the set
$\operatorname*{Inv}\sigma$ has at least $k-1$ elements. In other words,
$\left\vert \operatorname*{Inv}\sigma\right\vert \geq k-1$. Thus, $\ell\left(
\sigma\right)  =\left\vert \operatorname*{Inv}\sigma\right\vert \geq k-1$.
Since $\sigma=\operatorname*{cyc}\nolimits_{i_{1},i_{2},\ldots,i_{k}}$, this
rewrites as $\ell\left(  \operatorname*{cyc}\nolimits_{i_{1},i_{2}%
,\ldots,i_{k}}\right)  \geq k-1$. This solves Exercise \ref{exe.perm.cycles}
\textbf{(c)}.

\textbf{(d)} \textit{First solution to Exercise \ref{exe.perm.cycles}
\textbf{(d)}:} Let $i_{1},i_{2},\ldots,i_{k}$ be $k$ distinct elements of
$\left[  n\right]  $. Thus, $\left(  i_{1},i_{2},\ldots,i_{k}\right)  $ is a
list of some elements of $\left[  n\right]  $ such that $i_{1},i_{2}%
,\ldots,i_{k}$ are distinct. Hence, Proposition \ref{prop.perms.lists}
\textbf{(c)} (applied to $\left(  p_{1},p_{2},\ldots,p_{k}\right)  =\left(
i_{1},i_{2},\ldots,i_{k}\right)  $) yields that there exists a permutation
$\sigma\in S_{n}$ such that $\left(  i_{1},i_{2},\ldots,i_{k}\right)  =\left(
\sigma\left(  1\right)  ,\sigma\left(  2\right)  ,\ldots,\sigma\left(
k\right)  \right)  $. Consider such a $\sigma$.

Exercise \ref{exe.perm.cycles} \textbf{(b)} yields $\ell\left(
\operatorname*{cyc}\nolimits_{1,2,\ldots,k}\right)  =k-1$. But the definition
of $\left(  -1\right)  ^{\operatorname*{cyc}\nolimits_{1,2,\ldots,k}}$ yields
$\left(  -1\right)  ^{\operatorname*{cyc}\nolimits_{1,2,\ldots,k}}=\left(
-1\right)  ^{\ell\left(  \operatorname*{cyc}\nolimits_{1,2,\ldots,k}\right)
}=\left(  -1\right)  ^{k-1}$ (since $\ell\left(  \operatorname*{cyc}%
\nolimits_{1,2,\ldots,k}\right)  =k-1$).

Exercise \ref{exe.perm.cycles} \textbf{(a)} (applied to $1,2,\ldots,k$ instead
of $i_{1},i_{2},\ldots,i_{k}$) yields%
\[
\sigma\circ\operatorname*{cyc}\nolimits_{1,2,\ldots,k}\circ\sigma
^{-1}=\operatorname*{cyc}\nolimits_{\sigma\left(  1\right)  ,\sigma\left(
2\right)  ,\ldots,\sigma\left(  k\right)  }=\operatorname*{cyc}%
\nolimits_{i_{1},i_{2},\ldots,i_{k}}%
\]
(since $\left(  \sigma\left(  1\right)  ,\sigma\left(  2\right)
,\ldots,\sigma\left(  k\right)  \right)  =\left(  i_{1},i_{2},\ldots
,i_{k}\right)  $). Hence,
\[
\underbrace{\operatorname*{cyc}\nolimits_{i_{1},i_{2},\ldots,i_{k}}}%
_{=\sigma\circ\operatorname*{cyc}\nolimits_{1,2,\ldots,k}\circ\sigma^{-1}%
}\circ\sigma=\sigma\circ\operatorname*{cyc}\nolimits_{1,2,\ldots,k}%
\circ\underbrace{\sigma^{-1}\circ\sigma}_{=\operatorname*{id}}=\sigma
\circ\operatorname*{cyc}\nolimits_{1,2,\ldots,k}.
\]
Thus,%
\begin{align*}
\left(  -1\right)  ^{\operatorname*{cyc}\nolimits_{i_{1},i_{2},\ldots,i_{k}%
}\circ\sigma}  &  =\left(  -1\right)  ^{\sigma\circ\operatorname*{cyc}%
\nolimits_{1,2,\ldots,k}}=\left(  -1\right)  ^{\sigma}\cdot\underbrace{\left(
-1\right)  ^{\operatorname*{cyc}\nolimits_{1,2,\ldots,k}}}_{=\left(
-1\right)  ^{k-1}}\\
&  \ \ \ \ \ \ \ \ \ \ \left(  \text{by (\ref{eq.sign.prod}), applied to }%
\tau=\operatorname*{cyc}\nolimits_{1,2,\ldots,k}\right) \\
&  =\left(  -1\right)  ^{\sigma}\cdot\left(  -1\right)  ^{k-1}.
\end{align*}
Compared with%
\begin{align*}
\left(  -1\right)  ^{\operatorname*{cyc}\nolimits_{i_{1},i_{2},\ldots,i_{k}%
}\circ\sigma}  &  =\left(  -1\right)  ^{\operatorname*{cyc}\nolimits_{i_{1}%
,i_{2},\ldots,i_{k}}}\cdot\left(  -1\right)  ^{\sigma}\\
&  \ \ \ \ \ \ \ \ \ \ \left(  \text{by (\ref{eq.sign.prod}), applied to
}\operatorname*{cyc}\nolimits_{i_{1},i_{2},\ldots,i_{k}}\text{ and }%
\sigma\text{ instead of }\sigma\text{ and }\tau\right) \\
&  =\left(  -1\right)  ^{\sigma}\cdot\left(  -1\right)  ^{\operatorname*{cyc}%
\nolimits_{i_{1},i_{2},\ldots,i_{k}}},
\end{align*}
this yields $\left(  -1\right)  ^{\sigma}\cdot\left(  -1\right)
^{\operatorname*{cyc}\nolimits_{i_{1},i_{2},\ldots,i_{k}}}=\left(  -1\right)
^{\sigma}\cdot\left(  -1\right)  ^{k-1}$. We can cancel $\left(  -1\right)
^{\sigma}$ from this equality (since $\left(  -1\right)  ^{\sigma}\in\left\{
1,-1\right\}  $ is a nonzero integer), and thus obtain $\left(  -1\right)
^{\operatorname*{cyc}\nolimits_{i_{1},i_{2},\ldots,i_{k}}}=\left(  -1\right)
^{k-1}$. This solves Exercise \ref{exe.perm.cycles} \textbf{(d)}.
\end{verlong}

\textit{Second solution to Exercise \ref{exe.perm.cycles} \textbf{(d)}:} The
following solution to Exercise \ref{exe.perm.cycles} \textbf{(d)} appears in
various textbooks on abstract algebra.

\begin{verlong}
We have $k\in\left\{  1,2,\ldots,n\right\}  $, thus $k\geq1$, hence $k-1\geq
0$. Therefore, $k-1\in\mathbb{N}$.
\end{verlong}

Let $i_{1},i_{2},\ldots,i_{k}$ be $k$ distinct elements of $\left[  n\right]
$. We have%
\begin{equation}
\left(  -1\right)  ^{t_{i_{j},i_{j+1}}}=-1\ \ \ \ \ \ \ \ \ \ \text{for each
}j\in\left\{  1,2,\ldots,k-1\right\}  . \label{sol.perm.cycles.d.2nd.sign-t}%
\end{equation}

\begin{vershort}
\noindent Indeed, this follows from Exercise \ref{exe.ps4.1ab} \textbf{(b)}
(applied to $i_{j}$ and $i_{j+1}$ instead of $i$ and $j$).
\end{vershort}

\begin{verlong}
[\textit{Proof of (\ref{sol.perm.cycles.d.2nd.sign-t}):} Let $j\in\left\{
1,2,\ldots,k-1\right\}  $. Thus, $j+1\in\left\{  2,3,\ldots,k\right\}
\subseteq\left\{  1,2,\ldots,k\right\}  $ and $j\in\left\{  1,2,\ldots
,k-1\right\}  \subseteq\left\{  1,2,\ldots,k\right\}  $. Hence, $j$ and $j+1$
are two elements of $\left\{  1,2,\ldots,k\right\}  $. Thus, $i_{j}$ and
$i_{j+1}$ are two elements of $\left[  n\right]  $ (since $i_{1},i_{2}%
,\ldots,i_{k}$ are $k$ elements of $\left[  n\right]  $). Moreover, $j\neq
j+1$ and thus $i_{j}\neq i_{j+1}$ (since the $k$ elements $i_{1},i_{2}%
,\ldots,i_{k}$ are distinct). In other words, $i_{j}$ and $i_{j+1}$ are
distinct. Hence, $i_{j}$ and $i_{j+1}$ are two distinct elements of $\left[
n\right]  $. In other words, $i_{j}$ and $i_{j+1}$ are two distinct elements
of $\left\{  1,2,\ldots,n\right\}  $ (since $\left[  n\right]  =\left\{
1,2,\ldots,n\right\}  $). Thus, Exercise \ref{exe.ps4.1ab} \textbf{(b)}
(applied to $i_{j}$ and $i_{j+1}$ instead of $i$ and $j$) yields $\left(
-1\right)  ^{t_{i_{j},i_{j+1}}}=-1$. This proves
(\ref{sol.perm.cycles.d.2nd.sign-t}).]
\end{verlong}

Exercise \ref{exe.perm.c=ttt} yields%
\[
\operatorname*{cyc}\nolimits_{i_{1},i_{2},\ldots,i_{k}}=t_{i_{1},i_{2}}\circ
t_{i_{2},i_{3}}\circ\cdots\circ t_{i_{k-1},i_{k}}.
\]
Hence,%
\begin{align*}
\left(  -1\right)  ^{\operatorname*{cyc}\nolimits_{i_{1},i_{2},\ldots,i_{k}}}
&  =\left(  -1\right)  ^{t_{i_{1},i_{2}}\circ t_{i_{2},i_{3}}\circ\cdots\circ
t_{i_{k-1},i_{k}}}=\left(  -1\right)  ^{t_{i_{1},i_{2}}}\cdot\left(
-1\right)  ^{t_{i_{2},i_{3}}}\cdot\cdots\cdot\left(  -1\right)  ^{t_{i_{k-1}%
,i_{k}}}\\
&  \ \ \ \ \ \ \ \ \ \ \left(
\begin{array}
[c]{c}%
\text{by Proposition \ref{prop.sign.prod-of-many} (applied to }k-1\\
\text{and }t_{i_{j},i_{j+1}}\text{ instead of }k\text{ and }\sigma_{j}\text{)}%
\end{array}
\right) \\
&  =\prod_{j=1}^{k-1}\underbrace{\left(  -1\right)  ^{t_{i_{j},i_{j+1}}}%
}_{\substack{=-1\\\text{(by (\ref{sol.perm.cycles.d.2nd.sign-t}))}}%
}=\prod_{j=1}^{k-1}\left(  -1\right)  =\left(  -1\right)  ^{k-1}.
\end{align*}
This solves Exercise \ref{exe.perm.cycles} \textbf{(d)} again.
\end{proof}

\subsection{\label{sect.sol.perm.lehmer.prove}Solution to Exercise
\ref{exe.perm.lehmer.prove}}

We shall now approach the solution of Exercise \ref{exe.perm.lehmer.prove}.

Throughout Section \ref{sect.sol.perm.lehmer.prove}, we shall use the same
notations that were used in Section \ref{sect.perm.lehmer}. In particular, $n$
will denote a fixed element of $\mathbb{N}$.

We begin by proving Proposition \ref{prop.perm.lehmer.l}:

\begin{vershort}
\begin{proof}
[Proof of Proposition \ref{prop.perm.lehmer.l}.]Recall that the inversions of
$\sigma$ are defined to be the pairs $\left(  i,j\right)  $ of integers
satisfying $1\leq i<j\leq n$ and $\sigma\left(  i\right)  >\sigma\left(
j\right)  $. This definition can be rewritten as follows: The inversions of
$\sigma$ are the pairs $\left(  i,j\right)  \in\left[  n\right]  ^{2}$
satisfying $i<j$ and $\sigma\left(  i\right)  >\sigma\left(  j\right)  $.
Thus,%
\begin{align*}
&  \left(  \text{the number of inversions of }\sigma\right) \\
&  =\left(  \text{the number of all pairs }\left(  i,j\right)  \in\left[
n\right]  ^{2}\text{ satisfying }i<j\text{ and }\sigma\left(  i\right)
>\sigma\left(  j\right)  \right) \\
&  =\sum_{i\in\left[  n\right]  }\underbrace{\left(  \text{the number of all
}j\in\left[  n\right]  \text{ satisfying }i<j\text{ and }\sigma\left(
i\right)  >\sigma\left(  j\right)  \right)  }_{\substack{=\left(  \text{the
number of all }j\in\left\{  i+1,i+2,\ldots,n\right\}  \text{ and }%
\sigma\left(  i\right)  >\sigma\left(  j\right)  \right)  \\\text{(because the
}j\in\left[  n\right]  \text{ satisfying }i<j\text{ are precisely the }%
j\in\left\{  i+1,i+2,\ldots,n\right\}  \text{)}}}\\
&  =\sum_{i\in\left[  n\right]  }\underbrace{\left(  \text{the number of all
}j\in\left\{  i+1,i+2,\ldots,n\right\}  \text{ and }\sigma\left(  i\right)
>\sigma\left(  j\right)  \right)  }_{\substack{=\ell_{i}\left(  \sigma\right)
\\\text{(by the definition of }\ell_{i}\left(  \sigma\right)  \text{)}}%
}=\sum_{i\in\left[  n\right]  }\ell_{i}\left(  \sigma\right) \\
&  =\ell_{1}\left(  \sigma\right)  +\ell_{2}\left(  \sigma\right)
+\cdots+\ell_{n}\left(  \sigma\right)  .
\end{align*}
This proves Proposition \ref{prop.perm.lehmer.l}.
\end{proof}
\end{vershort}

\begin{verlong}
\begin{proof}
[Proof of Proposition \ref{prop.perm.lehmer.l}.]Let us use the Iverson bracket
notation (defined in Definition \ref{def.iverson}).

We have $\left[  n\right]  =\left\{  1,2,\ldots,n\right\}  $ (by the
definition of $\left[  n\right]  $).

We have%
\begin{equation}
\left\{  \text{inversions of }\sigma\right\}  =\left\{  \left(  u,v\right)
\in\left[  n\right]  ^{2}\ \mid\ u<v\text{ and }\sigma\left(  u\right)
>\sigma\left(  v\right)  \right\}  . \label{pf.prop.perm.lehmer.l.invsi=}%
\end{equation}

[\textit{Proof of (\ref{pf.prop.perm.lehmer.l.invsi=}):} Let $c\in\left\{
\text{inversions of }\sigma\right\}  $. We shall prove that \newline%
$c\in\left\{  \left(  u,v\right)  \in\left[  n\right]  ^{2}\ \mid\ u<v\text{
and }\sigma\left(  u\right)  >\sigma\left(  v\right)  \right\}  $.

We have $c\in\left\{  \text{inversions of }\sigma\right\}  $. In other words,
$c$ is an inversion of $\sigma$. In other words, $c$ is a pair $\left(
i,j\right)  $ of integers satisfying $1\leq i<j\leq n$ and $\sigma\left(
i\right)  >\sigma\left(  j\right)  $ (by the definition of an inversion).
Consider this $\left(  i,j\right)  $. Thus, $c=\left(  i,j\right)  $.

Both $i$ and $j$ are integers (since $\left(  i,j\right)  $ is a pair of
integers). We have $1\leq i\leq n$ and thus $i\in\left\{  1,2,\ldots
,n\right\}  $ (since $i$ is an integer). Hence, $i\in\left\{  1,2,\ldots
,n\right\}  =\left[  n\right]  $. Also, $1\leq j\leq n$ and thus $j\in\left\{
1,2,\ldots,n\right\}  $ (since $j$ is an integer). Hence, $j\in\left\{
1,2,\ldots,n\right\}  =\left[  n\right]  $. Combining $i\in\left[  n\right]  $
with $j\in\left[  n\right]  $, we obtain $\left(  i,j\right)  \in\left[
n\right]  \times\left[  n\right]  =\left[  n\right]  ^{2}$.

Now, $\left(  i,j\right)  $ is an element of $\left[  n\right]  ^{2}$ and
satisfies $i<j$ and $\sigma\left(  i\right)  >\sigma\left(  j\right)  $. In
other words, $\left(  i,j\right)  $ is a $\left(  u,v\right)  \in\left[
n\right]  ^{2}$ satisfying $u<v$ and $\sigma\left(  u\right)  >\sigma\left(
v\right)  $. In other words,%
\[
\left(  i,j\right)  \in\left\{  \left(  u,v\right)  \in\left[  n\right]
^{2}\ \mid\ u<v\text{ and }\sigma\left(  u\right)  >\sigma\left(  v\right)
\right\}  .
\]
Hence,%
\[
c=\left(  i,j\right)  \in\left\{  \left(  u,v\right)  \in\left[  n\right]
^{2}\ \mid\ u<v\text{ and }\sigma\left(  u\right)  >\sigma\left(  v\right)
\right\}  .
\]

Now, forget that we fixed $c$. We thus have proven that \newline$c\in\left\{
\left(  u,v\right)  \in\left[  n\right]  ^{2}\ \mid\ u<v\text{ and }%
\sigma\left(  u\right)  >\sigma\left(  v\right)  \right\}  $ for each
$c\in\left\{  \text{inversions of }\sigma\right\}  $. In other words,%
\begin{equation}
\left\{  \text{inversions of }\sigma\right\}  \subseteq\left\{  \left(
u,v\right)  \in\left[  n\right]  ^{2}\ \mid\ u<v\text{ and }\sigma\left(
u\right)  >\sigma\left(  v\right)  \right\}  .
\label{pf.prop.perm.lehmer.l.invsi=.pf.1}%
\end{equation}

On the other hand, let $d\in\left\{  \left(  u,v\right)  \in\left[  n\right]
^{2}\ \mid\ u<v\text{ and }\sigma\left(  u\right)  >\sigma\left(  v\right)
\right\}  $. We shall prove that $d\in\left\{  \text{inversions of }%
\sigma\right\}  $.

We have $d\in\left\{  \left(  u,v\right)  \in\left[  n\right]  ^{2}%
\ \mid\ u<v\text{ and }\sigma\left(  u\right)  >\sigma\left(  v\right)
\right\}  $. In other words, $d$ is a $\left(  u,v\right)  \in\left[
n\right]  ^{2}$ satisfying $u<v$ and $\sigma\left(  u\right)  >\sigma\left(
v\right)  $. Consider this $\left(  u,v\right)  $. Thus, $d=\left(
u,v\right)  $.

From $\left(  u,v\right)  \in\left[  n\right]  ^{2}$, we obtain $u\in\left[
n\right]  $ and $v\in\left[  n\right]  $. Thus, $u\in\left[  n\right]
=\left\{  1,2,\ldots,n\right\}  $, so that $1\leq u\leq n$. Also, $v\in\left[
n\right]  =\left\{  1,2,\ldots,n\right\}  $, so that $1\leq v\leq n$. Also,
$u$ is an integer (since $u\in\left\{  1,2,\ldots,n\right\}  $) and $v$ is an
integer (since $v\in\left\{  1,2,\ldots,n\right\}  $). Thus, $\left(
u,v\right)  $ is a pair of integers and satisfies $1\leq u<v\leq n$ and
$\sigma\left(  u\right)  >\sigma\left(  v\right)  $. In other words, $\left(
u,v\right)  $ is a pair $\left(  i,j\right)  $ of integers satisfying $1\leq
i<j\leq n$ and $\sigma\left(  i\right)  >\sigma\left(  j\right)  $. In other
words, $\left(  u,v\right)  $ is an inversion of $\sigma$ (by the definition
of an inversion). In other words, $\left(  u,v\right)  \in\left\{
\text{inversions of }\sigma\right\}  $. Thus, $d=\left(  u,v\right)
\in\left\{  \text{inversions of }\sigma\right\}  $.

Now, forget that we fixed $d$. We thus have proven that $d\in\left\{
\text{inversions of }\sigma\right\}  $ for each $d\in\left\{  \left(
u,v\right)  \in\left[  n\right]  ^{2}\ \mid\ u<v\text{ and }\sigma\left(
u\right)  >\sigma\left(  v\right)  \right\}  $. In other words,%
\[
\left\{  \left(  u,v\right)  \in\left[  n\right]  ^{2}\ \mid\ u<v\text{ and
}\sigma\left(  u\right)  >\sigma\left(  v\right)  \right\}  \subseteq\left\{
\text{inversions of }\sigma\right\}  .
\]
Combining this with (\ref{pf.prop.perm.lehmer.l.invsi=.pf.1}), we obtain%
\[
\left\{  \text{inversions of }\sigma\right\}  =\left\{  \left(  u,v\right)
\in\left[  n\right]  ^{2}\ \mid\ u<v\text{ and }\sigma\left(  u\right)
>\sigma\left(  v\right)  \right\}  .
\]
This proves (\ref{pf.prop.perm.lehmer.l.invsi=}).]

The equality (\ref{pf.prop.perm.lehmer.l.invsi=}) becomes%
\begin{align}
\left\{  \text{inversions of }\sigma\right\}   &  =\left\{  \left(
u,v\right)  \in\left[  n\right]  ^{2}\ \mid\ u<v\text{ and }\sigma\left(
u\right)  >\sigma\left(  v\right)  \right\} \nonumber\\
&  =\left\{  \left(  i,j\right)  \in\left[  n\right]  ^{2}\ \mid\ i<j\text{
and }\sigma\left(  i\right)  >\sigma\left(  j\right)  \right\}
\label{pf.prop.perm.lehmer.l.invsi=2}%
\end{align}
(here, we have renamed the index $\left(  u,v\right)  $ as $\left(
i,j\right)  $).

Recall that $\ell\left(  \sigma\right)  $ is defined to be the number of
inversions of $\sigma$. Thus,%
\begin{align}
\ell\left(  \sigma\right)   &  =\left(  \text{the number of inversions of
}\sigma\right) \nonumber\\
&  =\left\vert \underbrace{\left\{  \text{inversions of }\sigma\right\}
}_{\substack{=\left\{  \left(  i,j\right)  \in\left[  n\right]  ^{2}%
\ \mid\ i<j\text{ and }\sigma\left(  i\right)  >\sigma\left(  j\right)
\right\}  \\\text{(by (\ref{pf.prop.perm.lehmer.l.invsi=2}))}}}\right\vert
\nonumber\\
&  =\left\vert \left\{  \left(  i,j\right)  \in\left[  n\right]  ^{2}%
\ \mid\ i<j\text{ and }\sigma\left(  i\right)  >\sigma\left(  j\right)
\right\}  \right\vert . \label{pf.prop.perm.lehmer.l.l=}%
\end{align}

Now,%
\begin{align}
&  \underbrace{\sum_{i\in\left[  n\right]  }\sum_{j\in\left[  n\right]  }%
}_{=\sum_{\left(  i,j\right)  \in\left[  n\right]  ^{2}}}\left[  i<j\text{ and
}\sigma\left(  i\right)  >\sigma\left(  j\right)  \right] \nonumber\\
&  =\sum_{\left(  i,j\right)  \in\left[  n\right]  ^{2}}\left[  i<j\text{ and
}\sigma\left(  i\right)  >\sigma\left(  j\right)  \right] \nonumber\\
&  =\sum_{\substack{\left(  i,j\right)  \in\left[  n\right]  ^{2};\\i<j\text{
and }\sigma\left(  i\right)  >\sigma\left(  j\right)  }}\underbrace{\left[
i<j\text{ and }\sigma\left(  i\right)  >\sigma\left(  j\right)  \right]
}_{\substack{=1\\\text{(since we have }\left(  i<j\text{ and }\sigma\left(
i\right)  >\sigma\left(  j\right)  \right)  \text{)}}}\nonumber\\
&  \ \ \ \ \ \ \ \ \ \ +\sum_{\substack{\left(  i,j\right)  \in\left[
n\right]  ^{2};\\\text{not }\left(  i<j\text{ and }\sigma\left(  i\right)
>\sigma\left(  j\right)  \right)  }}\underbrace{\left[  i<j\text{ and }%
\sigma\left(  i\right)  >\sigma\left(  j\right)  \right]  }%
_{\substack{=0\\\text{(since we do not have }\left(  i<j\text{ and }%
\sigma\left(  i\right)  >\sigma\left(  j\right)  \right)  \text{)}%
}}\nonumber\\
&  \ \ \ \ \ \ \ \ \ \ \left(
\begin{array}
[c]{c}%
\text{since each }\left(  i,j\right)  \in\left[  n\right]  ^{2}\text{
satisfies either }\left(  i<j\text{ and }\sigma\left(  i\right)
>\sigma\left(  j\right)  \right) \\
\text{or }\left(  \text{not }\left(  i<j\text{ and }\sigma\left(  i\right)
>\sigma\left(  j\right)  \right)  \right)  \text{ (but not both)}%
\end{array}
\right) \nonumber\\
&  =\sum_{\substack{\left(  i,j\right)  \in\left[  n\right]  ^{2};\\i<j\text{
and }\sigma\left(  i\right)  >\sigma\left(  j\right)  }}1+\underbrace{\sum
_{\substack{\left(  i,j\right)  \in\left[  n\right]  ^{2};\\\text{not }\left(
i<j\text{ and }\sigma\left(  i\right)  >\sigma\left(  j\right)  \right)  }%
}0}_{=0}=\sum_{\substack{\left(  i,j\right)  \in\left[  n\right]
^{2};\\i<j\text{ and }\sigma\left(  i\right)  >\sigma\left(  j\right)
}}1\nonumber\\
&  =\left\vert \left\{  \left(  i,j\right)  \in\left[  n\right]  ^{2}%
\ \mid\ i<j\text{ and }\sigma\left(  i\right)  >\sigma\left(  j\right)
\right\}  \right\vert \cdot1\nonumber\\
&  =\left\vert \left\{  \left(  i,j\right)  \in\left[  n\right]  ^{2}%
\ \mid\ i<j\text{ and }\sigma\left(  i\right)  >\sigma\left(  j\right)
\right\}  \right\vert =\ell\left(  \sigma\right)
\label{pf.prop.perm.lehmer.l.3}%
\end{align}
(by (\ref{pf.prop.perm.lehmer.l.l=})).

On the other hand, each $i\in\left[  n\right]  $ satisfies%
\begin{equation}
\sum_{j\in\left[  n\right]  }\left[  i<j\text{ and }\sigma\left(  i\right)
>\sigma\left(  j\right)  \right]  =\ell_{i}\left(  \sigma\right)  .
\label{pf.prop.perm.lehmer.l.li=}%
\end{equation}

[\textit{Proof of (\ref{pf.prop.perm.lehmer.l.li=}):} Let $i\in\left[
n\right]  $. Recall that $\ell_{i}\left(  \sigma\right)  $ is the number of
all $j\in\left\{  i+1,i+2,\ldots,n\right\}  $ such that $\sigma\left(
i\right)  >\sigma\left(  j\right)  $ (by the definition of $\ell_{i}\left(
\sigma\right)  $). Thus,%
\begin{align}
\ell_{i}\left(  \sigma\right)   &  =\left(  \text{the number of all }%
j\in\left\{  i+1,i+2,\ldots,n\right\}  \text{ such that }\sigma\left(
i\right)  >\sigma\left(  j\right)  \right) \nonumber\\
&  =\left\vert \left\{  j\in\left\{  i+1,i+2,\ldots,n\right\}  \ \mid
\ \sigma\left(  i\right)  >\sigma\left(  j\right)  \right\}  \right\vert .
\label{pf.prop.perm.lehmer.l.li=.pf.1}%
\end{align}

On the other hand, $i\in\left[  n\right]  =\left\{  1,2,\ldots,n\right\}  $.
Thus, $1\leq i\leq n$. Now,%
\begin{align}
&  \underbrace{\sum_{j\in\left[  n\right]  }}_{\substack{=\sum_{j\in\left\{
1,2,\ldots,n\right\}  }\\\text{(since }\left[  n\right]  =\left\{
1,2,\ldots,n\right\}  \text{)}}}\left[  i<j\text{ and }\sigma\left(  i\right)
>\sigma\left(  j\right)  \right] \nonumber\\
&  =\underbrace{\sum_{j\in\left\{  1,2,\ldots,n\right\}  }}_{=\sum_{j=1}^{n}%
}\left[  i<j\text{ and }\sigma\left(  i\right)  >\sigma\left(  j\right)
\right]  =\sum_{j=1}^{n}\left[  i<j\text{ and }\sigma\left(  i\right)
>\sigma\left(  j\right)  \right] \nonumber\\
&  =\sum_{j=1}^{i}\left[  i<j\text{ and }\sigma\left(  i\right)
>\sigma\left(  j\right)  \right]  +\sum_{j=i+1}^{n}\left[  i<j\text{ and
}\sigma\left(  i\right)  >\sigma\left(  j\right)  \right]
\label{pf.prop.perm.lehmer.l.li=.pf.splitsum}%
\end{align}
(here, we have split the sum at $j=i$ (because $1\leq i\leq n$)).

For each $j\in\left\{  i+1,i+2,\ldots,n\right\}  $, we have%
\begin{equation}
\left[  i<j\text{ and }\sigma\left(  i\right)  >\sigma\left(  j\right)
\right]  =\left[  \sigma\left(  i\right)  >\sigma\left(  j\right)  \right]
\label{pf.prop.perm.lehmer.l.li=.pf.2}%
\end{equation}
\footnote{\textit{Proof of (\ref{pf.prop.perm.lehmer.l.li=.pf.2}):} Let
$j\in\left\{  i+1,i+2,\ldots,n\right\}  $. Thus, $j\geq i+1>i$, so that $i<j$.
We are in one of the following two cases:
\par
\textit{Case 1:} We have $\sigma\left(  i\right)  >\sigma\left(  j\right)  $.
\par
\textit{Case 2:} We don't have $\sigma\left(  i\right)  >\sigma\left(
j\right)  $.
\par
Let us first consider Case 1. In this case, we have $\sigma\left(  i\right)
>\sigma\left(  j\right)  $. Thus, we have $\left(  i<j\text{ and }%
\sigma\left(  i\right)  >\sigma\left(  j\right)  \right)  $ (since we also
have $i<j$). Hence, $\left[  i<j\text{ and }\sigma\left(  i\right)
>\sigma\left(  j\right)  \right]  =1$. But $\left[  \sigma\left(  i\right)
>\sigma\left(  j\right)  \right]  =1$ (since we have $\sigma\left(  i\right)
>\sigma\left(  j\right)  $). Thus, $\left[  i<j\text{ and }\sigma\left(
i\right)  >\sigma\left(  j\right)  \right]  =1=\left[  \sigma\left(  i\right)
>\sigma\left(  j\right)  \right]  $. Hence,
(\ref{pf.prop.perm.lehmer.l.li=.pf.2}) holds in Case 1.
\par
Let us now consider Case 2. In this case, we don't have $\sigma\left(
i\right)  >\sigma\left(  j\right)  $. Thus, we don't have $\left(  i<j\text{
and }\sigma\left(  i\right)  >\sigma\left(  j\right)  \right)  $ either.
Hence, $\left[  i<j\text{ and }\sigma\left(  i\right)  >\sigma\left(
j\right)  \right]  =0$. But $\left[  \sigma\left(  i\right)  >\sigma\left(
j\right)  \right]  =0$ (since we don't have $\sigma\left(  i\right)
>\sigma\left(  j\right)  $). Thus, $\left[  i<j\text{ and }\sigma\left(
i\right)  >\sigma\left(  j\right)  \right]  =0=\left[  \sigma\left(  i\right)
>\sigma\left(  j\right)  \right]  $. Hence,
(\ref{pf.prop.perm.lehmer.l.li=.pf.2}) holds in Case 2.
\par
We have now proven that (\ref{pf.prop.perm.lehmer.l.li=.pf.2}) holds in both
Cases 1 and 2. Since these two Cases cover all possibilities, we thus conclude
that (\ref{pf.prop.perm.lehmer.l.li=.pf.2}) always holds. Qed.}.

For each $j\in\left\{  1,2,\ldots,i\right\}  $, we have%
\begin{equation}
\left[  i<j\text{ and }\sigma\left(  i\right)  >\sigma\left(  j\right)
\right]  =0 \label{pf.prop.perm.lehmer.l.li=.pf.3}%
\end{equation}
\footnote{\textit{Proof of (\ref{pf.prop.perm.lehmer.l.li=.pf.3}):} Let
$j\in\left\{  1,2,\ldots,i\right\}  $. Thus, $j\leq i$. Hence, we don't have
$i<j$. Thus, we don't have $\left(  i<j\text{ and }\sigma\left(  i\right)
>\sigma\left(  j\right)  \right)  $ either. Hence, $\left[  i<j\text{ and
}\sigma\left(  i\right)  >\sigma\left(  j\right)  \right]  =0$. This proves
(\ref{pf.prop.perm.lehmer.l.li=.pf.3}).}.

Now, (\ref{pf.prop.perm.lehmer.l.li=.pf.splitsum}) becomes%
\begin{align*}
&  \sum_{j\in\left[  n\right]  }\left[  i<j\text{ and }\sigma\left(  i\right)
>\sigma\left(  j\right)  \right] \\
&  =\sum_{j=1}^{i}\underbrace{\left[  i<j\text{ and }\sigma\left(  i\right)
>\sigma\left(  j\right)  \right]  }_{\substack{=0\\\text{(by
(\ref{pf.prop.perm.lehmer.l.li=.pf.3}))}}}+\underbrace{\sum_{j=i+1}^{n}%
}_{=\sum_{j\in\left\{  i+1,i+2,\ldots,n\right\}  }}\underbrace{\left[
i<j\text{ and }\sigma\left(  i\right)  >\sigma\left(  j\right)  \right]
}_{\substack{=\left[  \sigma\left(  i\right)  >\sigma\left(  j\right)
\right]  \\\text{(by (\ref{pf.prop.perm.lehmer.l.li=.pf.2}))}}}\\
&  =\underbrace{\sum_{j=1}^{i}0}_{=0}+\sum_{j\in\left\{  i+1,i+2,\ldots
,n\right\}  }\left[  \sigma\left(  i\right)  >\sigma\left(  j\right)  \right]
=\sum_{j\in\left\{  i+1,i+2,\ldots,n\right\}  }\left[  \sigma\left(  i\right)
>\sigma\left(  j\right)  \right] \\
&  =\sum_{\substack{j\in\left\{  i+1,i+2,\ldots,n\right\}  ;\\\sigma\left(
i\right)  >\sigma\left(  j\right)  }}\underbrace{\left[  \sigma\left(
i\right)  >\sigma\left(  j\right)  \right]  }_{\substack{=1\\\text{(since
}\sigma\left(  i\right)  >\sigma\left(  j\right)  \text{)}}}+\sum
_{\substack{j\in\left\{  i+1,i+2,\ldots,n\right\}  ;\\\text{not }\sigma\left(
i\right)  >\sigma\left(  j\right)  }}\underbrace{\left[  \sigma\left(
i\right)  >\sigma\left(  j\right)  \right]  }_{\substack{=0\\\text{(since we
don't have }\sigma\left(  i\right)  >\sigma\left(  j\right)  \text{)}}}\\
&  \ \ \ \ \ \ \ \ \ \ \left(
\begin{array}
[c]{c}%
\text{since each }j\in\left\{  i+1,i+2,\ldots,n\right\}  \text{ satisfies
either }\left(  \sigma\left(  i\right)  >\sigma\left(  j\right)  \right) \\
\text{or }\left(  \text{not }\sigma\left(  i\right)  >\sigma\left(  j\right)
\right)  \text{ (but not both)}%
\end{array}
\right) \\
&  =\sum_{\substack{j\in\left\{  i+1,i+2,\ldots,n\right\}  ;\\\sigma\left(
i\right)  >\sigma\left(  j\right)  }}1+\underbrace{\sum_{\substack{j\in
\left\{  i+1,i+2,\ldots,n\right\}  ;\\\text{not }\sigma\left(  i\right)
>\sigma\left(  j\right)  }}0}_{=0}=\sum_{\substack{j\in\left\{  i+1,i+2,\ldots
,n\right\}  ;\\\sigma\left(  i\right)  >\sigma\left(  j\right)  }}1\\
&  =\left\vert \left\{  j\in\left\{  i+1,i+2,\ldots,n\right\}  \ \mid
\ \sigma\left(  i\right)  >\sigma\left(  j\right)  \right\}  \right\vert
\cdot1\\
&  =\left\vert \left\{  j\in\left\{  i+1,i+2,\ldots,n\right\}  \ \mid
\ \sigma\left(  i\right)  >\sigma\left(  j\right)  \right\}  \right\vert
=\ell_{i}\left(  \sigma\right)
\end{align*}
(by (\ref{pf.prop.perm.lehmer.l.li=.pf.1})). This proves
(\ref{pf.prop.perm.lehmer.l.li=}).]

Now, (\ref{pf.prop.perm.lehmer.l.3}) yields%
\begin{align*}
&  \ell\left(  \sigma\right) \\
&  =\sum_{i\in\left[  n\right]  }\underbrace{\sum_{j\in\left[  n\right]
}\left[  i<j\text{ and }\sigma\left(  i\right)  >\sigma\left(  j\right)
\right]  }_{\substack{=\ell_{i}\left(  \sigma\right)  \\\text{(by
(\ref{pf.prop.perm.lehmer.l.li=}))}}}=\underbrace{\sum_{i\in\left[  n\right]
}}_{\substack{=\sum_{i\in\left\{  1,2,\ldots,n\right\}  }\\\text{(since
}\left[  n\right]  =\left\{  1,2,\ldots,n\right\}  \text{)}}}\ell_{i}\left(
\sigma\right) \\
&  =\underbrace{\sum_{i\in\left\{  1,2,\ldots,n\right\}  }}_{=\sum_{i=1}^{n}%
}\ell_{i}\left(  \sigma\right)  =\sum_{i=1}^{n}\ell_{i}\left(  \sigma\right)
=\ell_{1}\left(  \sigma\right)  +\ell_{2}\left(  \sigma\right)  +\cdots
+\ell_{n}\left(  \sigma\right)  .
\end{align*}
This proves Proposition \ref{prop.perm.lehmer.l}.
\end{proof}
\end{verlong}

\begin{proof}
[Proof of Proposition \ref{prop.perm.lehmer.wd}.]Recall that $H=\left[
n-1\right]  _{0}\times\left[  n-2\right]  _{0}\times\cdots\times\left[
n-n\right]  _{0}$.

For each $i\in\left\{  1,2,\ldots,n\right\}  $, we have
\begin{equation}
\ell_{i}\left(  \sigma\right)  \in\left[  n-i\right]  _{0}.
\label{pf.prop.perm.lehmer.wd.1}%
\end{equation}

\begin{vershort}
[\textit{Proof of (\ref{pf.prop.perm.lehmer.wd.1}):} Let $i\in\left\{
1,2,\ldots,n\right\}  $. Thus, $i\in\left\{  1,2,\ldots,n\right\}  =\left[
n\right]  $. Recall that $\ell_{i}\left(  \sigma\right)  $ is the number of
all $j\in\left\{  i+1,i+2,\ldots,n\right\}  $ such that $\sigma\left(
i\right)  >\sigma\left(  j\right)  $ (by the definition of $\ell_{i}\left(
\sigma\right)  $). Thus,%
\begin{align*}
\ell_{i}\left(  \sigma\right)   &  =\left\vert \underbrace{\left\{
j\in\left\{  i+1,i+2,\ldots,n\right\}  \ \mid\ \sigma\left(  i\right)
>\sigma\left(  j\right)  \right\}  }_{\subseteq\left\{  i+1,i+2,\ldots
,n\right\}  }\right\vert \\
&  \leq\left\vert \left\{  i+1,i+2,\ldots,n\right\}  \right\vert
=n-i\ \ \ \ \ \ \ \ \ \ \left(  \text{since }i\leq n\right)  .
\end{align*}
Hence, $\ell_{i}\left(  \sigma\right)  \in\left\{  0,1,\ldots,n-i\right\}  $
(since $\ell_{i}\left(  \sigma\right)  \in\mathbb{N}$). But the definition of
$\left[  n-i\right]  _{0}$ yields $\left[  n-i\right]  _{0}=\left\{
0,1,\ldots,n-i\right\}  $. Hence, $\ell_{i}\left(  \sigma\right)  \in\left\{
0,1,\ldots,n-i\right\}  =\left[  n-i\right]  _{0}$. This proves
(\ref{pf.prop.perm.lehmer.wd.1}).]
\end{vershort}

\begin{verlong}
[\textit{Proof of (\ref{pf.prop.perm.lehmer.wd.1}):} Let $i\in\left\{
1,2,\ldots,n\right\}  $. Thus, $i\leq n$. Recall that $\left[  n\right]
=\left\{  1,2,\ldots,n\right\}  $ (by the definition of $\left[  n\right]  $).
Thus, $i\in\left\{  1,2,\ldots,n\right\}  =\left[  n\right]  $.

Recall that $\ell_{i}\left(  \sigma\right)  $ is the number of all
$j\in\left\{  i+1,i+2,\ldots,n\right\}  $ such that $\sigma\left(  i\right)
>\sigma\left(  j\right)  $ (by the definition of $\ell_{i}\left(
\sigma\right)  $). Thus,%
\begin{align*}
\ell_{i}\left(  \sigma\right)   &  =\left(  \text{the number of all }%
j\in\left\{  i+1,i+2,\ldots,n\right\}  \text{ such that }\sigma\left(
i\right)  >\sigma\left(  j\right)  \right) \\
&  =\left\vert \underbrace{\left\{  j\in\left\{  i+1,i+2,\ldots,n\right\}
\ \mid\ \sigma\left(  i\right)  >\sigma\left(  j\right)  \right\}
}_{\subseteq\left\{  i+1,i+2,\ldots,n\right\}  }\right\vert \\
&  \leq\left\vert \left\{  i+1,i+2,\ldots,n\right\}  \right\vert
=n-i\ \ \ \ \ \ \ \ \ \ \left(  \text{since }i\leq n\right)  .
\end{align*}
On the other hand, $\ell_{i}\left(  \sigma\right)  =\left\vert \left\{
j\in\left\{  i+1,i+2,\ldots,n\right\}  \ \mid\ \sigma\left(  i\right)
>\sigma\left(  j\right)  \right\}  \right\vert \in\mathbb{N}$. Hence, from
$\ell_{i}\left(  \sigma\right)  \leq n-i$, we obtain $\ell_{i}\left(
\sigma\right)  \in\left\{  0,1,\ldots,n-i\right\}  $.

But the definition of $\left[  n-i\right]  _{0}$ yields $\left[  n-i\right]
_{0}=\left\{  0,1,\ldots,n-i\right\}  $. Hence, $\ell_{i}\left(
\sigma\right)  \in\left\{  0,1,\ldots,n-i\right\}  =\left[  n-i\right]  _{0}$.
This proves (\ref{pf.prop.perm.lehmer.wd.1}).]
\end{verlong}

From (\ref{pf.prop.perm.lehmer.wd.1}), we know that $\ell_{i}\left(
\sigma\right)  \in\left[  n-i\right]  _{0}$ for each $i\in\left\{
1,2,\ldots,n\right\}  $. In other words,%
\[
\left(  \ell_{1}\left(  \sigma\right)  ,\ell_{2}\left(  \sigma\right)
,\ldots,\ell_{n}\left(  \sigma\right)  \right)  \in\left[  n-1\right]
_{0}\times\left[  n-2\right]  _{0}\times\cdots\times\left[  n-n\right]  _{0}.
\]
In view of $H=\left[  n-1\right]  _{0}\times\left[  n-2\right]  _{0}%
\times\cdots\times\left[  n-n\right]  _{0}$, this rewrites as \newline$\left(
\ell_{1}\left(  \sigma\right)  ,\ell_{2}\left(  \sigma\right)  ,\ldots
,\ell_{n}\left(  \sigma\right)  \right)  \in H$. This proves Proposition
\ref{prop.perm.lehmer.wd}.
\end{proof}

\begin{vershort}
\begin{proof}
[Proof of Lemma \ref{lem.perm.lexico1.lis}.]\textbf{(a)} We know that
$\ell_{i}\left(  \sigma\right)  $ is the number of all $j\in\left\{
i+1,i+2,\ldots,n\right\}  $ such that $\sigma\left(  i\right)  >\sigma\left(
j\right)  $ (by the definition of $\ell_{i}\left(  \sigma\right)  $). Hence,%
\begin{align}
\ell_{i}\left(  \sigma\right)   &  =\left(  \text{the number of all }%
j\in\left\{  i+1,i+2,\ldots,n\right\}  \text{ such that }\sigma\left(
i\right)  >\sigma\left(  j\right)  \right) \nonumber\\
&  =\left\vert \left\{  j\in\left\{  i+1,i+2,\ldots,n\right\}  \ \mid
\ \sigma\left(  i\right)  >\sigma\left(  j\right)  \right\}  \right\vert .
\label{pf.lem.perm.lexico1.lis.1}%
\end{align}
Define a set $A$ by
\begin{equation}
A=\left\{  j\in\left\{  i+1,i+2,\ldots,n\right\}  \ \mid\ \sigma\left(
i\right)  >\sigma\left(  j\right)  \right\}  .
\label{pf.lem.perm.lexico1.lis.A=}%
\end{equation}
Thus,%
\begin{equation}
\left\vert A\right\vert =\left\vert \left\{  j\in\left\{  i+1,i+2,\ldots
,n\right\}  \ \mid\ \sigma\left(  i\right)  >\sigma\left(  j\right)  \right\}
\right\vert =\ell_{i}\left(  \sigma\right)
\label{pf.lem.perm.lexico1.lis.sizeA=}%
\end{equation}
(by (\ref{pf.lem.perm.lexico1.lis.1})).

Let $B$ be the set $\left[  \sigma\left(  i\right)  -1\right]  \setminus
\sigma\left(  \left[  i\right]  \right)  $.

The map $\sigma$ is a permutation of $\left[  n\right]  $ (since $\sigma\in
S_{n}$), and thus is invertible, and therefore is injective.

For each $k\in A$, we have $\sigma\left(  k\right)  \in B$%
\ \ \ \ \footnote{\textit{Proof.} Let $k\in A$. Thus, $k\in A=\left\{
j\in\left\{  i+1,i+2,\ldots,n\right\}  \ \mid\ \sigma\left(  i\right)
>\sigma\left(  j\right)  \right\}  $. In other words, $k$ is an element of
$\left\{  i+1,i+2,\ldots,n\right\}  $ and satisfies $\sigma\left(  i\right)
>\sigma\left(  k\right)  $.
\par
From $k\in\left\{  i+1,i+2,\ldots,n\right\}  \subseteq\left[  n\right]  $, we
conclude that $\sigma\left(  k\right)  $ is well-defined. Also, $\sigma\left(
k\right)  <\sigma\left(  i\right)  $ (since $\sigma\left(  i\right)
>\sigma\left(  k\right)  $), so that $\sigma\left(  k\right)  \leq
\sigma\left(  i\right)  -1$ (since $\sigma\left(  k\right)  $ and
$\sigma\left(  i\right)  $ are integers). Thus, $\sigma\left(  k\right)
\in\left[  \sigma\left(  i\right)  -1\right]  $.
\par
Next, let us prove that $\sigma\left(  k\right)  \notin\sigma\left(  \left[
i\right]  \right)  $.
\par
Indeed, assume the contrary (for the sake of contradiction). Hence,
$\sigma\left(  k\right)  \in\sigma\left(  \left[  i\right]  \right)  $. In
other words, there exists some $j\in\left[  i\right]  $ such that
$\sigma\left(  k\right)  =\sigma\left(  j\right)  $. Consider this $j$. From
$\sigma\left(  k\right)  =\sigma\left(  j\right)  $, we obtain $k=j$ (since
the map $\sigma$ is injective). Hence, $k=j\in\left[  i\right]  $. But
$k\in\left\{  i+1,i+2,\ldots,n\right\}  =\left[  n\right]  \setminus\left[
i\right]  $, so that $k\notin\left[  i\right]  $. This contradicts
$k\in\left[  i\right]  $. This contradiction shows that our assumption was
false. Hence, $\sigma\left(  k\right)  \notin\sigma\left(  \left[  i\right]
\right)  $ is proven.
\par
Combining $\sigma\left(  k\right)  \in\left[  \sigma\left(  i\right)
-1\right]  $ with $\sigma\left(  k\right)  \notin\sigma\left(  \left[
i\right]  \right)  $, we obtain $\sigma\left(  k\right)  \in\left[
\sigma\left(  i\right)  -1\right]  \setminus\sigma\left(  \left[  i\right]
\right)  =B$. Qed.}. Hence, we can define a map $\alpha:A\rightarrow B$ by
\[
\left(  \alpha\left(  k\right)  =\sigma\left(  k\right)
\ \ \ \ \ \ \ \ \ \ \text{for each }k\in A\right)  .
\]
Consider this $\alpha$.

On the other hand, for each $k\in B$, we have $\sigma^{-1}\left(  k\right)
\in A$\ \ \ \ \footnote{\textit{Proof.} Let $k\in B$. Thus, $k\in B=\left[
\sigma\left(  i\right)  -1\right]  \setminus\sigma\left(  \left[  i\right]
\right)  $. In other words, $k\in\left[  \sigma\left(  i\right)  -1\right]  $
and $k\notin\sigma\left(  \left[  i\right]  \right)  $.
\par
From $k\in\left[  \sigma\left(  i\right)  -1\right]  $, we obtain $1\leq
k\leq\sigma\left(  i\right)  -1$. Also, $k\in\left[  \sigma\left(  i\right)
-1\right]  \subseteq\left[  n\right]  $, so that $\sigma^{-1}\left(  k\right)
$ is a well-defined element of $\left[  n\right]  $.
\par
We have $\sigma\left(  \sigma^{-1}\left(  k\right)  \right)  =k\leq
\sigma\left(  i\right)  -1<\sigma\left(  i\right)  $. In other words,
$\sigma\left(  i\right)  >\sigma\left(  \sigma^{-1}\left(  k\right)  \right)
$.
\par
Next, we claim that $\sigma^{-1}\left(  k\right)  \in\left\{  i+1,i+2,\ldots
,n\right\}  $. Indeed, assume the contrary (for the sake of contradiction).
Thus, $\sigma^{-1}\left(  k\right)  \notin\left\{  i+1,i+2,\ldots,n\right\}
$. Combining this with $\sigma^{-1}\left(  k\right)  \in\left[  n\right]  $,
we obtain%
\[
\sigma^{-1}\left(  k\right)  \in\left[  n\right]  \setminus\left\{
i+1,i+2,\ldots,n\right\}  =\left[  i\right]  .
\]
Hence, $k=\sigma\left(  \underbrace{\sigma^{-1}\left(  k\right)  }_{\in\left[
i\right]  }\right)  \in\sigma\left(  \left[  i\right]  \right)  $, which
contradicts $k\notin\sigma\left(  \left[  i\right]  \right)  $. This
contradiction shows that our assumption was false. Thus, $\sigma^{-1}\left(
k\right)  \in\left\{  i+1,i+2,\ldots,n\right\}  $ is proven.
\par
Now, we know that $\sigma^{-1}\left(  k\right)  \in\left\{  i+1,i+2,\ldots
,n\right\}  $ and $\sigma\left(  i\right)  >\sigma\left(  \sigma^{-1}\left(
k\right)  \right)  $. In other words, $\sigma^{-1}\left(  k\right)  $ is a
$j\in\left\{  i+1,i+2,\ldots,n\right\}  $ satisfying $\sigma\left(  i\right)
>\sigma\left(  j\right)  $. In other words,%
\[
\sigma^{-1}\left(  k\right)  \in\left\{  j\in\left\{  i+1,i+2,\ldots
,n\right\}  \ \mid\ \sigma\left(  i\right)  >\sigma\left(  j\right)  \right\}
.
\]
In view of (\ref{pf.lem.perm.lexico1.lis.A=}), this rewrites as $\sigma
^{-1}\left(  k\right)  \in A$. Qed.}. Hence, we can define a map
$\beta:B\rightarrow A$ by%
\[
\left(  \beta\left(  k\right)  =\sigma^{-1}\left(  k\right)
\ \ \ \ \ \ \ \ \ \ \text{for each }k\in B\right)  .
\]
Consider this $\beta$.

The maps $\alpha$ and $\beta$ are mutually inverse (since $\alpha$ is a
restriction of $\sigma$, whereas $\beta$ is a restriction of $\sigma^{-1}$),
and therefore are bijections. Hence, there is a bijection from $A$ to $B$
(namely, $\alpha$). Thus, $\left\vert A\right\vert =\left\vert B\right\vert $.

But (\ref{pf.lem.perm.lexico1.lis.sizeA=}) yields%
\[
\ell_{i}\left(  \sigma\right)  =\left\vert A\right\vert =\left\vert
B\right\vert =\left\vert \left[  \sigma\left(  i\right)  -1\right]
\setminus\sigma\left(  \left[  i\right]  \right)  \right\vert
\]
(since $B=\left[  \sigma\left(  i\right)  -1\right]  \setminus\sigma\left(
\left[  i\right]  \right)  $). This proves Lemma \ref{lem.perm.lexico1.lis}
\textbf{(a)}.

\textbf{(b)} If we had $\sigma\left(  i\right)  \in\left[  \sigma\left(
i\right)  -1\right]  $, then we would have $\sigma\left(  i\right)  \leq
\sigma\left(  i\right)  -1<\sigma\left(  i\right)  $, which would be absurd.
Hence, we have $\sigma\left(  i\right)  \notin\left[  \sigma\left(  i\right)
-1\right]  $.

But $\left[  i\right]  =\left\{  i\right\}  \cup\left[  i-1\right]  $. Hence,%
\[
\sigma\left(  \underbrace{\left[  i\right]  }_{=\left\{  i\right\}
\cup\left[  i-1\right]  }\right)  =\sigma\left(  \left\{  i\right\}
\cup\left[  i-1\right]  \right)  =\underbrace{\sigma\left(  \left\{
i\right\}  \right)  }_{=\left\{  \sigma\left(  i\right)  \right\}  }\cup
\sigma\left(  \left[  i-1\right]  \right)  =\left\{  \sigma\left(  i\right)
\right\}  \cup\sigma\left(  \left[  i-1\right]  \right)  .
\]
Thus,%
\begin{align*}
&  \left[  \sigma\left(  i\right)  -1\right]  \setminus\underbrace{\sigma
\left(  \left[  i\right]  \right)  }_{=\left\{  \sigma\left(  i\right)
\right\}  \cup\sigma\left(  \left[  i-1\right]  \right)  }\\
&  =\left[  \sigma\left(  i\right)  -1\right]  \setminus\left(  \left\{
\sigma\left(  i\right)  \right\}  \cup\sigma\left(  \left[  i-1\right]
\right)  \right) \\
&  =\underbrace{\left(  \left[  \sigma\left(  i\right)  -1\right]
\setminus\left\{  \sigma\left(  i\right)  \right\}  \right)  }%
_{\substack{=\left[  \sigma\left(  i\right)  -1\right]  \\\text{(since }%
\sigma\left(  i\right)  \notin\left[  \sigma\left(  i\right)  -1\right]
\text{)}}}\setminus\sigma\left(  \left[  i-1\right]  \right)  =\left[
\sigma\left(  i\right)  -1\right]  \setminus\sigma\left(  \left[  i-1\right]
\right)  .
\end{align*}

Now, Lemma \ref{lem.perm.lexico1.lis} \textbf{(a)} yields%
\[
\ell_{i}\left(  \sigma\right)  =\left\vert \underbrace{\left[  \sigma\left(
i\right)  -1\right]  \setminus\sigma\left(  \left[  i\right]  \right)
}_{=\left[  \sigma\left(  i\right)  -1\right]  \setminus\sigma\left(  \left[
i-1\right]  \right)  }\right\vert =\left\vert \left[  \sigma\left(  i\right)
-1\right]  \setminus\sigma\left(  \left[  i-1\right]  \right)  \right\vert .
\]
This proves Lemma \ref{lem.perm.lexico1.lis} \textbf{(b)}.

\textbf{(c)} The map $\sigma$ is a permutation of $\left[  n\right]  $, and
thus injective. Thus, Lemma \ref{lem.jectivity.pigeon0} \textbf{(c)} (applied
to $U=\left[  n\right]  $, $V=\left[  n\right]  $, $f=\sigma$ and $S=\left[
i-1\right]  $) yields $\left\vert \sigma\left(  \left[  i-1\right]  \right)
\right\vert =\left\vert \left[  i-1\right]  \right\vert $.

Now, Lemma \ref{lem.perm.lexico1.lis} \textbf{(b)} yields%
\begin{align*}
\ell_{i}\left(  \sigma\right)   &  =\left\vert \left[  \sigma\left(  i\right)
-1\right]  \setminus\sigma\left(  \left[  i-1\right]  \right)  \right\vert
\geq\underbrace{\left\vert \left[  \sigma\left(  i\right)  -1\right]
\right\vert }_{\substack{=\sigma\left(  i\right)  -1\\\text{(since }\left\vert
\left[  k\right]  \right\vert =k\\\text{for every }k\in\mathbb{N}\text{)}%
}}-\underbrace{\left\vert \sigma\left(  \left[  i-1\right]  \right)
\right\vert }_{\substack{=\left\vert \left[  i-1\right]  \right\vert
=i-1\\\text{(since }\left\vert \left[  k\right]  \right\vert =k\\\text{for
every }k\in\mathbb{N}\text{)}}}\\
&  \ \ \ \ \ \ \ \ \ \ \left(  \text{since any two finite sets }A\text{ and
}B\text{ satisfies }\left\vert A\setminus B\right\vert \geq\left\vert
A\right\vert -\left\vert B\right\vert \right) \\
&  =\left(  \sigma\left(  i\right)  -1\right)  -\left(  i-1\right)
=\sigma\left(  i\right)  -i.
\end{align*}
In other words, $\sigma\left(  i\right)  \leq i+\ell_{i}\left(  \sigma\right)
$. This proves Lemma \ref{lem.perm.lexico1.lis} \textbf{(c)}.
\end{proof}
\end{vershort}

\begin{verlong}
\begin{proof}
[Proof of Lemma \ref{lem.perm.lexico1.lis}.]We have $i\in\left[  n\right]
=\left\{  1,2,\ldots,n\right\}  $ (by the definition of $\left[  n\right]  $).
Thus, $1\leq i\leq n$. Hence, $\left\{  i+1,i+2,\ldots,n\right\}
\subseteq\left\{  1,2,\ldots,n\right\}  =\left[  n\right]  $. Also, the
definition of $\left[  i\right]  $ yields $\left[  i\right]  =\left\{
1,2,\ldots,i\right\}  \subseteq\left\{  1,2,\ldots,n\right\}  $ (since $i\leq
n$), so that $\left[  i\right]  \subseteq\left\{  1,2,\ldots,n\right\}
=\left[  n\right]  $.

Also, $i$ is a positive integer (since $i\in\left\{  1,2,\ldots,n\right\}  $).

We have $\sigma\in S_{n}$. In other words, $\sigma$ is a permutation of
$\left\{  1,2,\ldots,n\right\}  $ (since $S_{n}$ is the set of all
permutations of $\left\{  1,2,\ldots,n\right\}  $). In other words, $\sigma$
is a permutation of $\left[  n\right]  $ (since $\left[  n\right]  =\left\{
1,2,\ldots,n\right\}  $). In other words, $\sigma$ is a bijection $\left[
n\right]  \rightarrow\left[  n\right]  $. Hence, this map $\sigma$ is
invertible. Thus, the map $\sigma$ is bijective, and therefore injective.

Recall that $\sigma$ is a map $\left[  n\right]  \rightarrow\left[  n\right]
$. Thus, $\sigma\left(  \left[  i\right]  \right)  $ is a well-defined subset
of $\left[  n\right]  $ (since $\left[  i\right]  \subseteq\left[  n\right]  $).

Also, the definition of $\left[  i-1\right]  $ yields $\left[  i-1\right]
=\left\{  1,2,\ldots,i-1\right\}  \subseteq\left\{  1,2,\ldots,n\right\}  $
(since $i-1<i\leq n$), so that $\left[  i-1\right]  \subseteq\left\{
1,2,\ldots,n\right\}  =\left[  n\right]  $. Hence, $\sigma\left(  \left[
i-1\right]  \right)  $ is a well-defined subset of $\left[  n\right]  $ (since
$\sigma$ is a map $\left[  n\right]  \rightarrow\left[  n\right]  $).

We have%
\[
\underbrace{\left[  n\right]  }_{=\left\{  1,2,\ldots,n\right\}  }%
\setminus\underbrace{\left[  i\right]  }_{=\left\{  1,2,\ldots,i\right\}
}=\left\{  1,2,\ldots,n\right\}  \setminus\left\{  1,2,\ldots,i\right\}
=\left\{  i+1,i+2,\ldots,n\right\}  .
\]

We have $i\in\left[  n\right]  $ and thus $\sigma\left(  i\right)  \in\left[
n\right]  $ (since $\sigma$ is a map $\left[  n\right]  \rightarrow\left[
n\right]  $), so that $\sigma\left(  i\right)  \in\left[  n\right]  =\left\{
1,2,\ldots,n\right\}  $. Thus, $\sigma\left(  i\right)  \leq n$.

\textbf{(a)} We know that $\ell_{i}\left(  \sigma\right)  $ is the number of
all $j\in\left\{  i+1,i+2,\ldots,n\right\}  $ such that $\sigma\left(
i\right)  >\sigma\left(  j\right)  $ (by the definition of $\ell_{i}\left(
\sigma\right)  $). Hence,%
\begin{align}
\ell_{i}\left(  \sigma\right)   &  =\left(  \text{the number of all }%
j\in\left\{  i+1,i+2,\ldots,n\right\}  \text{ such that }\sigma\left(
i\right)  >\sigma\left(  j\right)  \right) \nonumber\\
&  =\left\vert \left\{  j\in\left\{  i+1,i+2,\ldots,n\right\}  \ \mid
\ \sigma\left(  i\right)  >\sigma\left(  j\right)  \right\}  \right\vert .
\label{pf.lem.perm.lexico1.lis.long.1}%
\end{align}

Define a set $A$ by
\begin{equation}
A=\left\{  j\in\left\{  i+1,i+2,\ldots,n\right\}  \ \mid\ \sigma\left(
i\right)  >\sigma\left(  j\right)  \right\}  .
\label{pf.lem.perm.lexico1.lis.long.A=}%
\end{equation}
Thus,%
\begin{equation}
\left\vert A\right\vert =\left\vert \left\{  j\in\left\{  i+1,i+2,\ldots
,n\right\}  \ \mid\ \sigma\left(  i\right)  >\sigma\left(  j\right)  \right\}
\right\vert =\ell_{i}\left(  \sigma\right)
\label{pf.lem.perm.lexico1.lis.long.sizeA=}%
\end{equation}
(by (\ref{pf.lem.perm.lexico1.lis.long.1})).

Also,%
\begin{align*}
A  &  =\left\{  j\in\left\{  i+1,i+2,\ldots,n\right\}  \ \mid\ \sigma\left(
i\right)  >\sigma\left(  j\right)  \right\} \\
&  \subseteq\left\{  i+1,i+2,\ldots,n\right\}  =\left[  n\right]
\setminus\left[  i\right]  \subseteq\left[  n\right]  .
\end{align*}

Define a set $B$ by%
\[
B=\left[  \sigma\left(  i\right)  -1\right]  \setminus\sigma\left(  \left[
i\right]  \right)  .
\]
Thus,%
\begin{align*}
B  &  =\left[  \sigma\left(  i\right)  -1\right]  \setminus\sigma\left(
\left[  i\right]  \right)  \subseteq\left[  \sigma\left(  i\right)  -1\right]
\\
&  =\left\{  1,2,\ldots,\sigma\left(  i\right)  -1\right\}
\ \ \ \ \ \ \ \ \ \ \left(  \text{by the definition of }\left[  \sigma\left(
i\right)  -1\right]  \right) \\
&  \subseteq\left\{  1,2,\ldots,n\right\}  \ \ \ \ \ \ \ \ \ \ \left(
\text{since }\sigma\left(  i\right)  -1\leq n\text{ (because }\sigma\left(
i\right)  -1<\sigma\left(  i\right)  \leq n\text{)}\right) \\
&  =\left[  n\right]  .
\end{align*}

For each $k\in A$, we have $\sigma\left(  k\right)  \in B$%
\ \ \ \ \footnote{\textit{Proof.} Let $k\in A$. Thus, $k\in A=\left\{
j\in\left\{  i+1,i+2,\ldots,n\right\}  \ \mid\ \sigma\left(  i\right)
>\sigma\left(  j\right)  \right\}  $. In other words, $k$ is a $j\in\left\{
i+1,i+2,\ldots,n\right\}  $ satisfying $\sigma\left(  i\right)  >\sigma\left(
j\right)  $. In other words, $k$ is an element of $\left\{  i+1,i+2,\ldots
,n\right\}  $ and satisfies $\sigma\left(  i\right)  >\sigma\left(  k\right)
$.
\par
From $k\in\left\{  i+1,i+2,\ldots,n\right\}  \subseteq\left[  n\right]  $, we
conclude that the element $\sigma\left(  k\right)  $ of $\left[  n\right]  $
is well-defined (since $\sigma$ is a bijection $\left[  n\right]
\rightarrow\left[  n\right]  $). Also, $\sigma\left(  k\right)  <\sigma\left(
i\right)  $ (since $\sigma\left(  i\right)  >\sigma\left(  k\right)  $), so
that $\sigma\left(  k\right)  \leq\sigma\left(  i\right)  -1$ (since
$\sigma\left(  k\right)  $ and $\sigma\left(  i\right)  $ are integers). But
$\sigma\left(  k\right)  \in\left[  n\right]  =\left\{  1,2,\ldots,n\right\}
$, so that $\sigma\left(  k\right)  \geq1$. Combining this with $\sigma\left(
k\right)  \leq\sigma\left(  i\right)  -1$, we obtain $1\leq\sigma\left(
k\right)  \leq\sigma\left(  i\right)  -1$. Hence, $\sigma\left(  k\right)
\in\left\{  1,2,\ldots,\sigma\left(  i\right)  -1\right\}  $. But $\left[
\sigma\left(  i\right)  -1\right]  =\left\{  1,2,\ldots,\sigma\left(
i\right)  -1\right\}  $ (by the definition of $\left[  \sigma\left(  i\right)
-1\right]  $). Thus, $\sigma\left(  k\right)  \in\left\{  1,2,\ldots
,\sigma\left(  i\right)  -1\right\}  =\left[  \sigma\left(  i\right)
-1\right]  $.
\par
Next, let us prove that $\sigma\left(  k\right)  \notin\sigma\left(  \left[
i\right]  \right)  $.
\par
Indeed, assume the contrary (for the sake of contradiction). Hence,
$\sigma\left(  k\right)  \in\sigma\left(  \left[  i\right]  \right)  $. In
other words, there exists some $j\in\left[  i\right]  $ such that
$\sigma\left(  k\right)  =\sigma\left(  j\right)  $. Consider this $j$. From
$\sigma\left(  k\right)  =\sigma\left(  j\right)  $, we obtain $k=j$ (since
the map $\sigma$ is injective). Hence, $k=j\in\left[  i\right]  $. But
$k\in\left\{  i+1,i+2,\ldots,n\right\}  =\left[  n\right]  \setminus\left[
i\right]  $. In other words, $k\in\left[  n\right]  $ and $k\notin\left[
i\right]  $. But $k\notin\left[  i\right]  $ contradicts $k\in\left[
i\right]  $. This contradiction shows that our assumption was false. Hence,
$\sigma\left(  k\right)  \notin\sigma\left(  \left[  i\right]  \right)  $ is
proven.
\par
Combining $\sigma\left(  k\right)  \in\left[  \sigma\left(  i\right)
-1\right]  $ with $\sigma\left(  k\right)  \notin\sigma\left(  \left[
i\right]  \right)  $, we obtain $\sigma\left(  k\right)  \in\left[
\sigma\left(  i\right)  -1\right]  \setminus\sigma\left(  \left[  i\right]
\right)  =B$. Qed.}. Hence, we can define a map $\alpha:A\rightarrow B$ by
\[
\left(  \alpha\left(  k\right)  =\sigma\left(  k\right)
\ \ \ \ \ \ \ \ \ \ \text{for each }k\in A\right)  .
\]
Consider this $\alpha$.

On the other hand, for each $k\in B$, we have $\sigma^{-1}\left(  k\right)
\in A$\ \ \ \ \footnote{\textit{Proof.} Let $k\in B$. Thus, $k\in B=\left[
\sigma\left(  i\right)  -1\right]  \setminus\sigma\left(  \left[  i\right]
\right)  $. In other words, $k\in\left[  \sigma\left(  i\right)  -1\right]  $
and $k\notin\sigma\left(  \left[  i\right]  \right)  $.
\par
We have $k\in\left[  \sigma\left(  i\right)  -1\right]  =\left\{
1,2,\ldots,\sigma\left(  i\right)  -1\right\}  $ (by the definition of
$\left[  \sigma\left(  i\right)  -1\right]  $), and therefore $1\leq
k\leq\sigma\left(  i\right)  -1$.
\par
Also, $k\in B\subseteq\left[  n\right]  $. Therefore, $\sigma^{-1}\left(
k\right)  $ is a well-defined element of $\left[  n\right]  $ (since $\sigma$
is a bijection $\left[  n\right]  \rightarrow\left[  n\right]  $).
\par
We have $\sigma\left(  \sigma^{-1}\left(  k\right)  \right)  =k\leq
\sigma\left(  i\right)  -1<\sigma\left(  i\right)  $. In other words,
$\sigma\left(  i\right)  >\sigma\left(  \sigma^{-1}\left(  k\right)  \right)
$.
\par
Next, we claim that $\sigma^{-1}\left(  k\right)  \in\left\{  i+1,i+2,\ldots
,n\right\}  $. Indeed, assume the contrary (for the sake of contradiction).
Thus, $\sigma^{-1}\left(  k\right)  \notin\left\{  i+1,i+2,\ldots,n\right\}
$. Combining this with $\sigma^{-1}\left(  k\right)  \in\left[  n\right]  $,
we obtain%
\[
\sigma^{-1}\left(  k\right)  \in\left[  n\right]  \setminus
\underbrace{\left\{  i+1,i+2,\ldots,n\right\}  }_{=\left[  n\right]
\setminus\left[  i\right]  }=\left[  n\right]  \setminus\left(  \left[
n\right]  \setminus\left[  i\right]  \right)  \subseteq\left[  i\right]  .
\]
Hence, $k=\sigma\left(  \underbrace{\sigma^{-1}\left(  k\right)  }_{\in\left[
i\right]  }\right)  \in\sigma\left(  \left[  i\right]  \right)  $, which
contradicts $k\notin\sigma\left(  \left[  i\right]  \right)  $. This
contradiction shows that our assumption was false. Thus, $\sigma^{-1}\left(
k\right)  \in\left\{  i+1,i+2,\ldots,n\right\}  $ is proven.
\par
Now, we know that $\sigma^{-1}\left(  k\right)  $ is an element of $\left\{
i+1,i+2,\ldots,n\right\}  $ and satisfies $\sigma\left(  i\right)
>\sigma\left(  \sigma^{-1}\left(  k\right)  \right)  $. In other words,
$\sigma^{-1}\left(  k\right)  $ is a $j\in\left\{  i+1,i+2,\ldots,n\right\}  $
satisfying $\sigma\left(  i\right)  >\sigma\left(  j\right)  $. In other
words,%
\[
\sigma^{-1}\left(  k\right)  \in\left\{  j\in\left\{  i+1,i+2,\ldots
,n\right\}  \ \mid\ \sigma\left(  i\right)  >\sigma\left(  j\right)  \right\}
.
\]
In view of (\ref{pf.lem.perm.lexico1.lis.long.A=}), this rewrites as
$\sigma^{-1}\left(  k\right)  \in A$. Qed.}. Hence, we can define a map
$\beta:B\rightarrow A$ by%
\[
\left(  \beta\left(  k\right)  =\sigma^{-1}\left(  k\right)
\ \ \ \ \ \ \ \ \ \ \text{for each }k\in B\right)  .
\]
Consider this $\beta$.

We have $\alpha\circ\beta=\operatorname*{id}$\ \ \ \ \footnote{\textit{Proof.}
Let $k\in B$. Then, $\beta\left(  k\right)  =\sigma^{-1}\left(  k\right)  $
(by the definition of $\beta$). But the definition of $\alpha$ yields
$\alpha\left(  \beta\left(  k\right)  \right)  =\sigma\left(
\underbrace{\beta\left(  k\right)  }_{=\sigma^{-1}\left(  k\right)  }\right)
=\sigma\left(  \sigma^{-1}\left(  k\right)  \right)  =k$. Hence, $\left(
\alpha\circ\beta\right)  \left(  k\right)  =\alpha\left(  \beta\left(
k\right)  \right)  =k=\operatorname*{id}\left(  k\right)  $.
\par
Now, forget that we fixed $k$. We thus have proven that $\left(  \alpha
\circ\beta\right)  \left(  k\right)  =\operatorname*{id}\left(  k\right)  $
for each $k\in B$. In other words, $\alpha\circ\beta=\operatorname*{id}$.} and
$\beta\circ\alpha=\operatorname*{id}$\ \ \ \ \footnote{\textit{Proof.} Let
$k\in A$. Then, $\alpha\left(  k\right)  =\sigma\left(  k\right)  $ (by the
definition of $\alpha$). But the definition of $\beta$ yields $\beta\left(
\alpha\left(  k\right)  \right)  =\sigma^{-1}\left(  \underbrace{\alpha\left(
k\right)  }_{=\sigma\left(  k\right)  }\right)  =\sigma^{-1}\left(
\sigma\left(  k\right)  \right)  =k$. Hence, $\left(  \beta\circ\alpha\right)
\left(  k\right)  =\beta\left(  \alpha\left(  k\right)  \right)
=k=\operatorname*{id}\left(  k\right)  $.
\par
Now, forget that we fixed $k$. We thus have proven that $\left(  \beta
\circ\alpha\right)  \left(  k\right)  =\operatorname*{id}\left(  k\right)  $
for each $k\in A$. In other words, $\beta\circ\alpha=\operatorname*{id}$.}.
These two equalities show that the maps $\alpha$ and $\beta$ are mutually
inverse. Hence, the map $\alpha$ is invertible, i.e., is bijective. In other
words, the map $\alpha$ is a bijection. Hence, there is a bijection from $A$
to $B$ (namely, $\alpha$). Thus, $\left\vert A\right\vert =\left\vert
B\right\vert $.

But (\ref{pf.lem.perm.lexico1.lis.long.sizeA=}) yields%
\[
\ell_{i}\left(  \sigma\right)  =\left\vert A\right\vert =\left\vert
B\right\vert =\left\vert \left[  \sigma\left(  i\right)  -1\right]
\setminus\sigma\left(  \left[  i\right]  \right)  \right\vert
\]
(since $B=\left[  \sigma\left(  i\right)  -1\right]  \setminus\sigma\left(
\left[  i\right]  \right)  $). This proves Lemma \ref{lem.perm.lexico1.lis}
\textbf{(a)}.

\textbf{(b)} If we had $\sigma\left(  i\right)  \in\left\{  1,2,\ldots
,\sigma\left(  i\right)  -1\right\}  $, then we would have $\sigma\left(
i\right)  \leq\sigma\left(  i\right)  -1<\sigma\left(  i\right)  $, which
would be absurd. Hence, we cannot have $\sigma\left(  i\right)  \in\left\{
1,2,\ldots,\sigma\left(  i\right)  -1\right\}  $. In other words, we have
$\sigma\left(  i\right)  \notin\left\{  1,2,\ldots,\sigma\left(  i\right)
-1\right\}  $.

The definition of $\left[  \sigma\left(  i\right)  -1\right]  $ yields
$\left[  \sigma\left(  i\right)  -1\right]  =\left\{  1,2,\ldots,\sigma\left(
i\right)  -1\right\}  $. Thus, $\sigma\left(  i\right)  \notin\left\{
1,2,\ldots,\sigma\left(  i\right)  -1\right\}  =\left[  \sigma\left(
i\right)  -1\right]  $.

But
\begin{align*}
\left[  i\right]   &  =\left\{  1,2,\ldots,i\right\}  =\underbrace{\left\{
1,2,\ldots,i-1\right\}  }_{=\left[  i-1\right]  }\cup\left\{  i\right\}
\ \ \ \ \ \ \ \ \ \ \left(  \text{since }i\text{ is a positive integer}\right)
\\
&  =\left[  i-1\right]  \cup\left\{  i\right\}  =\left\{  i\right\}
\cup\left[  i-1\right]  .
\end{align*}
Hence,%
\[
\sigma\left(  \underbrace{\left[  i\right]  }_{=\left\{  i\right\}
\cup\left[  i-1\right]  }\right)  =\sigma\left(  \left\{  i\right\}
\cup\left[  i-1\right]  \right)  =\underbrace{\sigma\left(  \left\{
i\right\}  \right)  }_{=\left\{  \sigma\left(  i\right)  \right\}  }\cup
\sigma\left(  \left[  i-1\right]  \right)  =\left\{  \sigma\left(  i\right)
\right\}  \cup\sigma\left(  \left[  i-1\right]  \right)  .
\]
Thus,%
\begin{align*}
&  \left[  \sigma\left(  i\right)  -1\right]  \setminus\underbrace{\sigma
\left(  \left[  i\right]  \right)  }_{=\left\{  \sigma\left(  i\right)
\right\}  \cup\sigma\left(  \left[  i-1\right]  \right)  }\\
&  =\left[  \sigma\left(  i\right)  -1\right]  \setminus\left(  \left\{
\sigma\left(  i\right)  \right\}  \cup\sigma\left(  \left[  i-1\right]
\right)  \right) \\
&  =\underbrace{\left(  \left[  \sigma\left(  i\right)  -1\right]
\setminus\left\{  \sigma\left(  i\right)  \right\}  \right)  }%
_{\substack{=\left[  \sigma\left(  i\right)  -1\right]  \\\text{(since }%
\sigma\left(  i\right)  \notin\left[  \sigma\left(  i\right)  -1\right]
\text{)}}}\setminus\sigma\left(  \left[  i-1\right]  \right)  =\left[
\sigma\left(  i\right)  -1\right]  \setminus\sigma\left(  \left[  i-1\right]
\right)  .
\end{align*}

Now, Lemma \ref{lem.perm.lexico1.lis} \textbf{(a)} yields%
\[
\ell_{i}\left(  \sigma\right)  =\left\vert \underbrace{\left[  \sigma\left(
i\right)  -1\right]  \setminus\sigma\left(  \left[  i\right]  \right)
}_{=\left[  \sigma\left(  i\right)  -1\right]  \setminus\sigma\left(  \left[
i-1\right]  \right)  }\right\vert =\left\vert \left[  \sigma\left(  i\right)
-1\right]  \setminus\sigma\left(  \left[  i-1\right]  \right)  \right\vert .
\]
This proves Lemma \ref{lem.perm.lexico1.lis} \textbf{(b)}.

\textbf{(c)} The map $\sigma:\left[  n\right]  \rightarrow\left[  n\right]  $
is injective. The set $\left[  n\right]  $ is finite. Also, $\left[
i-1\right]  \subseteq\left[  n\right]  $. That is, $\left[  i-1\right]  $ is a
subset of $\left[  n\right]  $. Thus, Lemma \ref{lem.jectivity.pigeon0}
\textbf{(c)} (applied to $U=\left[  n\right]  $, $V=\left[  n\right]  $,
$f=\sigma$ and $S=\left[  i-1\right]  $) yields $\left\vert \sigma\left(
\left[  i-1\right]  \right)  \right\vert =\left\vert \left[  i-1\right]
\right\vert $.

We have $i-1\in\mathbb{N}$ (since $i$ is a positive integer). Also,
$\sigma\left(  i\right)  $ is a positive integer (since $\sigma\left(
i\right)  \in\left\{  1,2,\ldots,n\right\}  $); hence, $\sigma\left(
i\right)  -1\in\mathbb{N}$.

For each $k\in\mathbb{N}$, we have $\left[  k\right]  =\left\{  1,2,\ldots
,k\right\}  $ (by the definition of $\left[  k\right]  $) and thus $\left\vert
\left[  k\right]  \right\vert =\left\vert \left\{  1,2,\ldots,k\right\}
\right\vert =k$ (since $k\in\mathbb{N}$). Applying this to $k=i-1$, we
conclude that $\left\vert \left[  i-1\right]  \right\vert =i-1$ (since
$i-1\in\mathbb{N}$).

For each $k\in\mathbb{N}$, we have $\left\vert \left[  k\right]  \right\vert
=k$ (as we have just proved). Applying this to $k=\sigma\left(  i\right)  -1$,
we conclude that $\left\vert \left[  \sigma\left(  i\right)  -1\right]
\right\vert =\sigma\left(  i\right)  -1$ (since $\sigma\left(  i\right)
-1\in\mathbb{N}$).

If $A$ and $B$ are two finite sets, then $\left\vert A\setminus B\right\vert
\geq\left\vert A\right\vert -\left\vert B\right\vert $%
\ \ \ \ \footnote{\textit{Proof.} Let $A$ and $B$ be two finite sets. Then,
$A\cap B$ is a subset of $B$. Hence, $\left\vert A\cap B\right\vert
\leq\left\vert B\right\vert $.
\par
But $A\cap B$ also is a subset of $A$. Hence, $\left\vert A\setminus\left(
A\cap B\right)  \right\vert =\left\vert A\right\vert -\left\vert A\cap
B\right\vert $. In view of $A\setminus\left(  A\cap B\right)  =A\setminus B$,
this rewrites as $\left\vert A\setminus B\right\vert =\left\vert A\right\vert
-\underbrace{\left\vert A\cap B\right\vert }_{\leq\left\vert B\right\vert
}\geq\left\vert A\right\vert -\left\vert B\right\vert $, qed.}. Applying this
to $A=\left[  \sigma\left(  i\right)  -1\right]  $ and $B=\sigma\left(
\left[  i-1\right]  \right)  $, we obtain
\begin{align*}
\left\vert \left[  \sigma\left(  i\right)  -1\right]  \setminus\sigma\left(
\left[  i-1\right]  \right)  \right\vert  &  \geq\underbrace{\left\vert
\left[  \sigma\left(  i\right)  -1\right]  \right\vert }_{\substack{=\sigma
\left(  i\right)  -1}}-\underbrace{\left\vert \sigma\left(  \left[
i-1\right]  \right)  \right\vert }_{=\left\vert \left[  i-1\right]
\right\vert =i-1}\\
&  =\left(  \sigma\left(  i\right)  -1\right)  -\left(  i-1\right)
=\sigma\left(  i\right)  -i.
\end{align*}

Now, Lemma \ref{lem.perm.lexico1.lis} \textbf{(b)} yields $\ell_{i}\left(
\sigma\right)  =\left\vert \left[  \sigma\left(  i\right)  -1\right]
\setminus\sigma\left(  \left[  i-1\right]  \right)  \right\vert \geq
\sigma\left(  i\right)  -i$. In other words, $i+\ell_{i}\left(  \sigma\right)
\geq\sigma\left(  i\right)  $; thus, $\sigma\left(  i\right)  \leq i+\ell
_{i}\left(  \sigma\right)  $. This proves Lemma \ref{lem.perm.lexico1.lis}
\textbf{(c)}.
\end{proof}
\end{verlong}

Next, we state a further fact, which will be used in proving Proposition
\ref{prop.perm.lehmer.lex}:

\begin{lemma}
\label{lem.perm.lehmer.lex1}Let $\sigma\in S_{n}$ and $\tau\in S_{n}$. Let
$i\in\left[  n\right]  $. Assume that
\begin{equation}
\text{each }k\in\left[  i-1\right]  \text{ satisfies }\sigma\left(  k\right)
=\tau\left(  k\right)  . \label{eq.lem.perm.lehmer.lex1.ass}%
\end{equation}
Then:

\textbf{(a)} Each $k\in\left[  i\right]  $ satisfies $\sigma\left(  \left[
k-1\right]  \right)  =\tau\left(  \left[  k-1\right]  \right)  $.

\textbf{(b)} We have $\ell_{k}\left(  \sigma\right)  =\ell_{k}\left(
\tau\right)  $ for each $k\in\left[  i-1\right]  $.

\textbf{(c)} Assume furthermore that $\sigma\left(  i\right)  <\tau\left(
i\right)  $. Then, $\ell_{i}\left(  \sigma\right)  <\ell_{i}\left(
\tau\right)  $.
\end{lemma}

\begin{vershort}
\begin{proof}
[Proof of Lemma \ref{lem.perm.lehmer.lex1}.]\textbf{(a)} Let $k\in\left[
i\right]  $. Thus, $k\leq i$.

Let $j\in\left[  k-1\right]  $. Thus, $j\leq\underbrace{k}_{\leq i}-1\leq
i-1$, so that $j\in\left[  i-1\right]  $. Hence,
(\ref{eq.lem.perm.lehmer.lex1.ass}) (applied to $j$ instead of $k$) shows that
$\sigma\left(  j\right)  =\tau\left(  j\right)  $.

Now, forget that we fixed $j$. We thus have shown that $\sigma\left(
j\right)  =\tau\left(  j\right)  $ for each $j\in\left[  k-1\right]  $. In
other words,%
\[
\left(  \sigma\left(  1\right)  ,\sigma\left(  2\right)  ,\ldots,\sigma\left(
k-1\right)  \right)  =\left(  \tau\left(  1\right)  ,\tau\left(  2\right)
,\ldots,\tau\left(  k-1\right)  \right)  .
\]
Thus,%
\[
\left\{  \sigma\left(  1\right)  ,\sigma\left(  2\right)  ,\ldots
,\sigma\left(  k-1\right)  \right\}  =\left\{  \tau\left(  1\right)
,\tau\left(  2\right)  ,\ldots,\tau\left(  k-1\right)  \right\}  .
\]

Now,%
\begin{align*}
\sigma\left(  \underbrace{\left[  k-1\right]  }_{=\left\{  1,2,\ldots
,k-1\right\}  }\right)   &  =\sigma\left(  \left\{  1,2,\ldots,k-1\right\}
\right)  =\left\{  \sigma\left(  1\right)  ,\sigma\left(  2\right)
,\ldots,\sigma\left(  k-1\right)  \right\} \\
&  =\left\{  \tau\left(  1\right)  ,\tau\left(  2\right)  ,\ldots,\tau\left(
k-1\right)  \right\}  =\tau\left(  \underbrace{\left\{  1,2,\ldots
,k-1\right\}  }_{=\left[  k-1\right]  }\right)  =\tau\left(  \left[
k-1\right]  \right)  .
\end{align*}
This proves Lemma \ref{lem.perm.lehmer.lex1} \textbf{(a)}.

\textbf{(b)} Let $k\in\left[  i-1\right]  $. Then, Lemma
\ref{lem.perm.lexico1.lis} \textbf{(b)} (applied to $k$ instead of $i$) yields
$\ell_{k}\left(  \sigma\right)  =\left\vert \left[  \sigma\left(  k\right)
-1\right]  \setminus\sigma\left(  \left[  k-1\right]  \right)  \right\vert $.
The same argument (applied to $\tau$ instead of $\sigma$) yields $\ell
_{k}\left(  \tau\right)  =\left\vert \left[  \tau\left(  k\right)  -1\right]
\setminus\tau\left(  \left[  k-1\right]  \right)  \right\vert $.

But $k\in\left[  i-1\right]  \subseteq\left[  i\right]  $. Hence, Lemma
\ref{lem.perm.lehmer.lex1} \textbf{(a)} yields $\sigma\left(  \left[
k-1\right]  \right)  =\tau\left(  \left[  k-1\right]  \right)  $. Also,
(\ref{eq.lem.perm.lehmer.lex1.ass}) yields $\sigma\left(  k\right)
=\tau\left(  k\right)  $. Hence,%
\[
\ell_{k}\left(  \sigma\right)  =\left\vert \left[  \underbrace{\sigma\left(
k\right)  }_{=\tau\left(  k\right)  }-1\right]  \setminus\underbrace{\sigma
\left(  \left[  k-1\right]  \right)  }_{=\tau\left(  \left[  k-1\right]
\right)  }\right\vert =\left\vert \left[  \tau\left(  k\right)  -1\right]
\setminus\tau\left(  \left[  k-1\right]  \right)  \right\vert =\ell_{k}\left(
\tau\right)  .
\]
This proves Lemma \ref{lem.perm.lehmer.lex1} \textbf{(b)}.

\textbf{(c)} We have $i\in\left[  i\right]  $. Hence, Lemma
\ref{lem.perm.lehmer.lex1} \textbf{(a)} (applied to $k=i$) yields
$\sigma\left(  \left[  i-1\right]  \right)  =\tau\left(  \left[  i-1\right]
\right)  $.

Also, $\sigma\left(  i\right)  <\tau\left(  i\right)  $, so that
$\sigma\left(  i\right)  -1<\tau\left(  i\right)  -1$ and therefore $\left[
\sigma\left(  i\right)  -1\right]  \subseteq\left[  \tau\left(  i\right)
-1\right]  $.

From $\sigma\left(  i\right)  <\tau\left(  i\right)  $, we also obtain
$\sigma\left(  i\right)  \leq\tau\left(  i\right)  -1$ (since $\sigma\left(
i\right)  $ and $\tau\left(  i\right)  $ are integers), and thus
$\sigma\left(  i\right)  \in\left[  \tau\left(  i\right)  -1\right]  $.

Also, $\sigma\left(  i\right)  \notin\sigma\left(  \left[  i-1\right]
\right)  $. [\textit{Proof:} Assume the contrary. Thus, $\sigma\left(
i\right)  \in\sigma\left(  \left[  i-1\right]  \right)  $. In other words,
$\sigma\left(  i\right)  =\sigma\left(  j\right)  $ for some $j\in\left[
i-1\right]  $. Consider this $j$. From $\sigma\left(  i\right)  =\sigma\left(
j\right)  $, we obtain $i=j$ (since $\sigma$ is injective), so that
$i=j\in\left[  i-1\right]  $ and thus $i\leq i-1<i$. But this is absurd.
Hence, we found a contradiction, so that $\sigma\left(  i\right)  \notin%
\sigma\left(  \left[  i-1\right]  \right)  $ is proven.]

If we had $\sigma\left(  i\right)  \in\left[  \sigma\left(  i\right)
-1\right]  $, then we would have $\sigma\left(  i\right)  \leq\sigma\left(
i\right)  -1<\sigma\left(  i\right)  $, which is absurd. Hence, we have
$\sigma\left(  i\right)  \notin\left[  \sigma\left(  i\right)  -1\right]  $.
Thus, also $\sigma\left(  i\right)  \notin\left[  \sigma\left(  i\right)
-1\right]  \setminus\sigma\left(  \left[  i-1\right]  \right)  $.

Combining $\sigma\left(  i\right)  \in\left[  \tau\left(  i\right)  -1\right]
$ with $\sigma\left(  i\right)  \notin\sigma\left(  \left[  i-1\right]
\right)  $, we obtain $\sigma\left(  i\right)  \in\left[  \tau\left(
i\right)  -1\right]  \setminus\sigma\left(  \left[  i-1\right]  \right)  $.

Now,%
\begin{equation}
\underbrace{\left[  \sigma\left(  i\right)  -1\right]  }_{\subseteq\left[
\tau\left(  i\right)  -1\right]  }\setminus\sigma\left(  \left[  i-1\right]
\right)  \subseteq\left[  \tau\left(  i\right)  -1\right]  \setminus
\sigma\left(  \left[  i-1\right]  \right)  .
\label{pf.lem.perm.lehmer.lex1.short.c.2}%
\end{equation}

Moreover, the set $\left[  \tau\left(  i\right)  -1\right]  \setminus
\sigma\left(  \left[  i-1\right]  \right)  $ contains $\sigma\left(  i\right)
$ (since $\sigma\left(  i\right)  \in\left[  \tau\left(  i\right)  -1\right]
\setminus\sigma\left(  \left[  i-1\right]  \right)  $), but the set $\left[
\sigma\left(  i\right)  -1\right]  \setminus\sigma\left(  \left[  i-1\right]
\right)  $ does not (since $\sigma\left(  i\right)  \notin\left[
\sigma\left(  i\right)  -1\right]  \setminus\sigma\left(  \left[  i-1\right]
\right)  $). Thus, these two sets are distinct. In other words, $\left[
\sigma\left(  i\right)  -1\right]  \setminus\sigma\left(  \left[  i-1\right]
\right)  \neq\left[  \tau\left(  i\right)  -1\right]  \setminus\sigma\left(
\left[  i-1\right]  \right)  $. Combining this with
(\ref{pf.lem.perm.lehmer.lex1.short.c.2}), we conclude that $\left[
\sigma\left(  i\right)  -1\right]  \setminus\sigma\left(  \left[  i-1\right]
\right)  $ is a \textbf{proper} subset of $\left[  \tau\left(  i\right)
-1\right]  \setminus\sigma\left(  \left[  i-1\right]  \right)  $.

But recall the following fundamental fact: If $P$ is a finite set, and if $Q$
is a proper subset of $P$, then $\left\vert Q\right\vert <\left\vert
P\right\vert $. Applying this to $P=\left[  \tau\left(  i\right)  -1\right]
\setminus\sigma\left(  \left[  i-1\right]  \right)  $ and $Q=\left[
\sigma\left(  i\right)  -1\right]  \setminus\sigma\left(  \left[  i-1\right]
\right)  $, we conclude that%
\[
\left\vert \left[  \sigma\left(  i\right)  -1\right]  \setminus\sigma\left(
\left[  i-1\right]  \right)  \right\vert <\left\vert \left[  \tau\left(
i\right)  -1\right]  \setminus\sigma\left(  \left[  i-1\right]  \right)
\right\vert
\]
(since $\left[  \sigma\left(  i\right)  -1\right]  \setminus\sigma\left(
\left[  i-1\right]  \right)  $ is a \textbf{proper} subset of $\left[
\tau\left(  i\right)  -1\right]  \setminus\sigma\left(  \left[  i-1\right]
\right)  $).

But Lemma \ref{lem.perm.lexico1.lis} \textbf{(b)} yields $\ell_{i}\left(
\sigma\right)  =\left\vert \left[  \sigma\left(  i\right)  -1\right]
\setminus\sigma\left(  \left[  i-1\right]  \right)  \right\vert $. The same
argument (applied to $\tau$ instead of $\sigma$) yields $\ell_{i}\left(
\tau\right)  =\left\vert \left[  \tau\left(  i\right)  -1\right]
\setminus\tau\left(  \left[  i-1\right]  \right)  \right\vert $. Hence,%
\begin{align*}
\ell_{i}\left(  \sigma\right)   &  =\left\vert \left[  \sigma\left(  i\right)
-1\right]  \setminus\sigma\left(  \left[  i-1\right]  \right)  \right\vert \\
&  <\left\vert \left[  \tau\left(  i\right)  -1\right]  \setminus
\underbrace{\sigma\left(  \left[  i-1\right]  \right)  }_{=\tau\left(  \left[
i-1\right]  \right)  }\right\vert =\left\vert \left[  \tau\left(  i\right)
-1\right]  \setminus\tau\left(  \left[  i-1\right]  \right)  \right\vert
=\ell_{i}\left(  \tau\right)  .
\end{align*}
This proves Lemma \ref{lem.perm.lehmer.lex1} \textbf{(c)}.
\end{proof}
\end{vershort}

\begin{verlong}
\begin{proof}
[Proof of Lemma \ref{lem.perm.lehmer.lex1}.]We have $i\in\left[  n\right]
=\left\{  1,2,\ldots,n\right\}  $ (by the definition of $\left[  n\right]  $).
Thus, $1\leq i\leq n$. Hence, $\left\{  i+1,i+2,\ldots,n\right\}
\subseteq\left\{  1,2,\ldots,n\right\}  =\left[  n\right]  $. Also, the
definition of $\left[  i\right]  $ yields $\left[  i\right]  =\left\{
1,2,\ldots,i\right\}  \subseteq\left\{  1,2,\ldots,n\right\}  $ (since $i\leq
n$), so that $\left[  i\right]  \subseteq\left\{  1,2,\ldots,n\right\}
=\left[  n\right]  $.

Also, $i$ is a positive integer (since $i\in\left\{  1,2,\ldots,n\right\}  $).
Thus, $i\in\left\{  1,2,\ldots,i\right\}  =\left[  i\right]  $.

Also, the definition of $\left[  i-1\right]  $ yields $\left[  i-1\right]
=\left\{  1,2,\ldots,i-1\right\}  \subseteq\left\{  1,2,\ldots,i\right\}  $
(since $i-1<i$), so that $\left[  i-1\right]  \subseteq\left\{  1,2,\ldots
,i\right\}  =\left[  i\right]  \subseteq\left[  n\right]  $.

We have $\sigma\in S_{n}$. In other words, $\sigma$ is a permutation of
$\left\{  1,2,\ldots,n\right\}  $ (since $S_{n}$ is the set of all
permutations of $\left\{  1,2,\ldots,n\right\}  $). In other words, $\sigma$
is a permutation of $\left[  n\right]  $ (since $\left[  n\right]  =\left\{
1,2,\ldots,n\right\}  $). In other words, $\sigma$ is a bijection $\left[
n\right]  \rightarrow\left[  n\right]  $. Hence, this map $\sigma$ is
invertible. Thus, the map $\sigma$ is bijective, and therefore injective.

\textbf{(a)} Let $k\in\left[  i\right]  $. Thus, $k\in\left[  i\right]
=\left\{  1,2,\ldots,i\right\}  $, so that $k\leq i$. Also, $k\in\left\{
1,2,\ldots,i\right\}  $, so that $k$ is a positive integer. Thus,
$k-1\in\mathbb{N}$.

Let $j\in\left\{  1,2,\ldots,k-1\right\}  $. Hence, $j\leq\underbrace{k}_{\leq
i}-1\leq i-1$. Combining this with $j\geq1$ (which follows from $j\in\left\{
1,2,\ldots,k-1\right\}  $), we obtain $1\leq j\leq i-1$. Hence, $j\in\left\{
1,2,\ldots,i-1\right\}  =\left[  i-1\right]  $. Hence,
(\ref{eq.lem.perm.lehmer.lex1.ass}) (applied to $j$ instead of $k$) shows that
$\sigma\left(  j\right)  =\tau\left(  j\right)  $.

Now, forget that we fixed $j$. We thus have shown that $\sigma\left(
j\right)  =\tau\left(  j\right)  $ for each $j\in\left\{  1,2,\ldots
,k-1\right\}  $. In other words,%
\[
\left(  \sigma\left(  1\right)  ,\sigma\left(  2\right)  ,\ldots,\sigma\left(
k-1\right)  \right)  =\left(  \tau\left(  1\right)  ,\tau\left(  2\right)
,\ldots,\tau\left(  k-1\right)  \right)  .
\]
Thus,%
\[
\left\{  \sigma\left(  1\right)  ,\sigma\left(  2\right)  ,\ldots
,\sigma\left(  k-1\right)  \right\}  =\left\{  \tau\left(  1\right)
,\tau\left(  2\right)  ,\ldots,\tau\left(  k-1\right)  \right\}  .
\]

Now, $\left[  k-1\right]  =\left\{  1,2,\ldots,k-1\right\}  $ (by the
definition of $\left[  k-1\right]  $). Applying the map $\sigma$ to both sides
of this equality, we find%
\begin{align*}
\sigma\left(  \left[  k-1\right]  \right)   &  =\sigma\left(  \left\{
1,2,\ldots,k-1\right\}  \right)  =\left\{  \sigma\left(  1\right)
,\sigma\left(  2\right)  ,\ldots,\sigma\left(  k-1\right)  \right\} \\
&  =\left\{  \tau\left(  1\right)  ,\tau\left(  2\right)  ,\ldots,\tau\left(
k-1\right)  \right\}  =\tau\left(  \underbrace{\left\{  1,2,\ldots
,k-1\right\}  }_{=\left[  k-1\right]  }\right)  =\tau\left(  \left[
k-1\right]  \right)  .
\end{align*}
This proves Lemma \ref{lem.perm.lehmer.lex1} \textbf{(a)}.

\textbf{(b)} Let $k\in\left[  i-1\right]  $. Then, $k\in\left[  i-1\right]
\subseteq\left[  n\right]  $. Hence, Lemma \ref{lem.perm.lexico1.lis}
\textbf{(b)} (applied to $k$ instead of $i$) yields $\ell_{k}\left(
\sigma\right)  =\left\vert \left[  \sigma\left(  k\right)  -1\right]
\setminus\sigma\left(  \left[  k-1\right]  \right)  \right\vert $. The same
argument (applied to $\tau$ instead of $\sigma$) yields $\ell_{k}\left(
\tau\right)  =\left\vert \left[  \tau\left(  k\right)  -1\right]
\setminus\tau\left(  \left[  k-1\right]  \right)  \right\vert $.

But $k\in\left[  i-1\right]  \subseteq\left[  i\right]  $. Hence, Lemma
\ref{lem.perm.lehmer.lex1} \textbf{(a)} yields $\sigma\left(  \left[
k-1\right]  \right)  =\tau\left(  \left[  k-1\right]  \right)  $. Also,
(\ref{eq.lem.perm.lehmer.lex1.ass}) yields $\sigma\left(  k\right)
=\tau\left(  k\right)  $. Hence,%
\[
\ell_{k}\left(  \sigma\right)  =\left\vert \left[  \underbrace{\sigma\left(
k\right)  }_{=\tau\left(  k\right)  }-1\right]  \setminus\underbrace{\sigma
\left(  \left[  k-1\right]  \right)  }_{=\tau\left(  \left[  k-1\right]
\right)  }\right\vert =\left\vert \left[  \tau\left(  k\right)  -1\right]
\setminus\tau\left(  \left[  k-1\right]  \right)  \right\vert =\ell_{k}\left(
\tau\right)  .
\]
This proves Lemma \ref{lem.perm.lehmer.lex1} \textbf{(b)}.

\textbf{(c)} We have $i\in\left[  i\right]  $. Hence, Lemma
\ref{lem.perm.lehmer.lex1} \textbf{(a)} (applied to $k=i$) yields
$\sigma\left(  \left[  i-1\right]  \right)  =\tau\left(  \left[  i-1\right]
\right)  $.

Recall that $\sigma$ is a map $\left[  n\right]  \rightarrow\left[  n\right]
$. Thus, $\sigma\left(  i\right)  \in\left[  n\right]  =\left\{
1,2,\ldots,n\right\}  $. Hence, $\sigma\left(  i\right)  $ is an integer
satisfying $1\leq\sigma\left(  i\right)  \leq n$. The same argument (applied
to $\tau$ instead of $\sigma$) shows that $\tau\left(  i\right)  $ is an
integer satisfying $1\leq\tau\left(  i\right)  \leq n$.

We have $\left[  \sigma\left(  i\right)  -1\right]  =\left\{  1,2,\ldots
,\sigma\left(  i\right)  -1\right\}  $ (by the definition of $\left[
\sigma\left(  i\right)  -1\right]  $) and $\left[  \tau\left(  i\right)
-1\right]  =\left\{  1,2,\ldots,\tau\left(  i\right)  -1\right\}  $ (by the
definition of $\left[  \tau\left(  i\right)  -1\right]  $).

Also, $\sigma\left(  i\right)  <\tau\left(  i\right)  $, so that
$\underbrace{\sigma\left(  i\right)  }_{<\tau\left(  i\right)  }-1<\tau\left(
i\right)  -1$ and therefore $\left\{  1,2,\ldots,\sigma\left(  i\right)
-1\right\}  \subseteq\left\{  1,2,\ldots,\tau\left(  i\right)  -1\right\}  $.
Thus,%
\[
\left[  \sigma\left(  i\right)  -1\right]  =\left\{  1,2,\ldots,\sigma\left(
i\right)  -1\right\}  \subseteq\left\{  1,2,\ldots,\tau\left(  i\right)
-1\right\}  =\left[  \tau\left(  i\right)  -1\right]  .
\]

From $\sigma\left(  i\right)  <\tau\left(  i\right)  $, we also obtain
$\sigma\left(  i\right)  \leq\tau\left(  i\right)  -1$ (since $\sigma\left(
i\right)  $ and $\tau\left(  i\right)  $ are integers). Hence, $1\leq
\sigma\left(  i\right)  \leq\tau\left(  i\right)  -1$. Thus, $\sigma\left(
i\right)  \in\left\{  1,2,\ldots,\tau\left(  i\right)  -1\right\}  =\left[
\tau\left(  i\right)  -1\right]  $.

Also, $\sigma\left(  i\right)  \notin\sigma\left(  \left[  i-1\right]
\right)  $. [\textit{Proof:} Assume the contrary. Thus, $\sigma\left(
i\right)  \in\sigma\left(  \left[  i-1\right]  \right)  $. In other words,
$\sigma\left(  i\right)  =\sigma\left(  j\right)  $ for some $j\in\left[
i-1\right]  $. Consider this $j$. From $\sigma\left(  i\right)  =\sigma\left(
j\right)  $, we obtain $i=j$ (since $\sigma$ is injective), so that
$i=j\in\left[  i-1\right]  =\left\{  1,2,\ldots,i-1\right\}  $ and thus $i\leq
i-1<i$. But this is absurd. Hence, we found a contradiction. Thus, our
assumption was wrong. This completes the proof of $\sigma\left(  i\right)
\notin\sigma\left(  \left[  i-1\right]  \right)  $.]

Combining $\sigma\left(  i\right)  \in\left[  \tau\left(  i\right)  -1\right]
$ with $\sigma\left(  i\right)  \notin\sigma\left(  \left[  i-1\right]
\right)  $, we obtain $\sigma\left(  i\right)  \in\left[  \tau\left(
i\right)  -1\right]  \setminus\sigma\left(  \left[  i-1\right]  \right)  $.

Also,
\begin{equation}
\left[  \sigma\left(  i\right)  -1\right]  \setminus\sigma\left(  \left[
i-1\right]  \right)  \neq\left[  \tau\left(  i\right)  -1\right]
\setminus\sigma\left(  \left[  i-1\right]  \right)  .
\label{pf.lem.perm.lehmer.lex1.c.1}%
\end{equation}
[\textit{Proof of (\ref{pf.lem.perm.lehmer.lex1.c.1}):} Assume the contrary.
Thus, $\left[  \sigma\left(  i\right)  -1\right]  \setminus\sigma\left(
\left[  i-1\right]  \right)  =\left[  \tau\left(  i\right)  -1\right]
\setminus\sigma\left(  \left[  i-1\right]  \right)  $. Hence,
\begin{align*}
\sigma\left(  i\right)   &  \in\left[  \tau\left(  i\right)  -1\right]
\setminus\sigma\left(  \left[  i-1\right]  \right)  =\left[  \sigma\left(
i\right)  -1\right]  \setminus\sigma\left(  \left[  i-1\right]  \right) \\
&  \subseteq\left[  \sigma\left(  i\right)  -1\right]  =\left\{
1,2,\ldots,\sigma\left(  i\right)  -1\right\}  .
\end{align*}
Hence, $\sigma\left(  i\right)  \leq\sigma\left(  i\right)  -1<\sigma\left(
i\right)  $, which is absurd. Thus, we found a contradiction. Thus, our
assumption was wrong. This completes the proof of
(\ref{pf.lem.perm.lehmer.lex1.c.1}).]

Now,%
\[
\underbrace{\left[  \sigma\left(  i\right)  -1\right]  }_{\subseteq\left[
\tau\left(  i\right)  -1\right]  }\setminus\sigma\left(  \left[  i-1\right]
\right)  \subseteq\left[  \tau\left(  i\right)  -1\right]  \setminus
\sigma\left(  \left[  i-1\right]  \right)  .
\]
Combining this with (\ref{pf.lem.perm.lehmer.lex1.c.1}), we conclude that
$\left[  \sigma\left(  i\right)  -1\right]  \setminus\sigma\left(  \left[
i-1\right]  \right)  $ is a \textbf{proper} subset of $\left[  \tau\left(
i\right)  -1\right]  \setminus\sigma\left(  \left[  i-1\right]  \right)  $.

But recall the following fundamental fact: If $P$ is a finite set, and if $Q$
is a proper subset of $P$, then $\left\vert Q\right\vert <\left\vert
P\right\vert $. Applying this to $P=\left[  \tau\left(  i\right)  -1\right]
\setminus\sigma\left(  \left[  i-1\right]  \right)  $ and $Q=\left[
\sigma\left(  i\right)  -1\right]  \setminus\sigma\left(  \left[  i-1\right]
\right)  $, we conclude that%
\[
\left\vert \left[  \sigma\left(  i\right)  -1\right]  \setminus\sigma\left(
\left[  i-1\right]  \right)  \right\vert <\left\vert \left[  \tau\left(
i\right)  -1\right]  \setminus\sigma\left(  \left[  i-1\right]  \right)
\right\vert
\]
(since $\left[  \sigma\left(  i\right)  -1\right]  \setminus\sigma\left(
\left[  i-1\right]  \right)  $ is a \textbf{proper} subset of $\left[
\tau\left(  i\right)  -1\right]  \setminus\sigma\left(  \left[  i-1\right]
\right)  $).

But Lemma \ref{lem.perm.lexico1.lis} \textbf{(b)} yields $\ell_{i}\left(
\sigma\right)  =\left\vert \left[  \sigma\left(  i\right)  -1\right]
\setminus\sigma\left(  \left[  i-1\right]  \right)  \right\vert $. The same
argument (applied to $\tau$ instead of $\sigma$) yields $\ell_{i}\left(
\tau\right)  =\left\vert \left[  \tau\left(  i\right)  -1\right]
\setminus\tau\left(  \left[  i-1\right]  \right)  \right\vert $. Hence,%
\begin{align*}
\ell_{i}\left(  \sigma\right)   &  =\left\vert \left[  \sigma\left(  i\right)
-1\right]  \setminus\sigma\left(  \left[  i-1\right]  \right)  \right\vert \\
&  <\left\vert \left[  \tau\left(  i\right)  -1\right]  \setminus
\underbrace{\sigma\left(  \left[  i-1\right]  \right)  }_{=\tau\left(  \left[
i-1\right]  \right)  }\right\vert =\left\vert \left[  \tau\left(  i\right)
-1\right]  \setminus\tau\left(  \left[  i-1\right]  \right)  \right\vert
=\ell_{i}\left(  \tau\right)  .
\end{align*}
This proves Lemma \ref{lem.perm.lehmer.lex1} \textbf{(c)}.
\end{proof}
\end{verlong}

We are now ready to prove Proposition \ref{prop.perm.lehmer.lex}:

\begin{vershort}
\begin{proof}
[Proof of Proposition \ref{prop.perm.lehmer.lex}.]We have assumed that
\[
\left(  \sigma\left(  1\right)  ,\sigma\left(  2\right)  ,\ldots,\sigma\left(
n\right)  \right)  <_{\operatorname*{lex}}\left(  \tau\left(  1\right)
,\tau\left(  2\right)  ,\ldots,\tau\left(  n\right)  \right)  .
\]
According to Definition \ref{def.perm.lehmer.lex-ord}, this means the
following: There exists some $k\in\left[  n\right]  $ such that $\sigma\left(
k\right)  \neq\tau\left(  k\right)  $, and the \textbf{smallest} such $k$
satisfies $\sigma\left(  k\right)  <\tau\left(  k\right)  $.

Let $i$ be the smallest such $k$. Thus, $i$ is an element of $\left[
n\right]  $ such that $\sigma\left(  i\right)  <\tau\left(  i\right)  $, but
\[
\text{each }k\in\left[  i-1\right]  \text{ satisfies }\sigma\left(  k\right)
=\tau\left(  k\right)
\]
(since $i$ is the \textbf{smallest} $k\in\left[  n\right]  $ such that
$\sigma\left(  k\right)  \neq\tau\left(  k\right)  $).

Thus, Lemma \ref{lem.perm.lehmer.lex1} \textbf{(b)} shows that%
\begin{equation}
\text{each }k\in\left[  i-1\right]  \text{ satisfies }\ell_{k}\left(
\sigma\right)  =\ell_{k}\left(  \tau\right)  .
\label{pf.prop.perm.lehmer.lex.1}%
\end{equation}
Furthermore, Lemma \ref{lem.perm.lehmer.lex1} \textbf{(c)} shows that
$\ell_{i}\left(  \sigma\right)  <\ell_{i}\left(  \tau\right)  $ (since
$\sigma\left(  i\right)  <\tau\left(  i\right)  $). Thus, $\ell_{i}\left(
\sigma\right)  \neq\ell_{i}\left(  \tau\right)  $. In other words, $i$ is a
$k\in\left[  n\right]  $ such that $\ell_{k}\left(  \sigma\right)  \neq
\ell_{k}\left(  \tau\right)  $. Moreover, (\ref{pf.prop.perm.lehmer.lex.1})
shows that $i$ is the \textbf{smallest} such $k$. Thus, the smallest
$k\in\left[  n\right]  $ such that $\ell_{k}\left(  \sigma\right)  \neq
\ell_{k}\left(  \tau\right)  $ satisfies $\ell_{k}\left(  \sigma\right)
<\ell_{k}\left(  \tau\right)  $ (because this $k$ is $i$, and $i$ satisfies
$\ell_{i}\left(  \sigma\right)  <\ell_{i}\left(  \tau\right)  $).

Thus, we have shown that there exists some $k\in\left[  n\right]  $ such that
$\ell_{k}\left(  \sigma\right)  \neq\ell_{k}\left(  \tau\right)  $, and the
\textbf{smallest} such $k$ satisfies $\ell_{k}\left(  \sigma\right)  <\ell
_{k}\left(  \tau\right)  $. But this means precisely that%
\[
\left(  \ell_{1}\left(  \sigma\right)  ,\ell_{2}\left(  \sigma\right)
,\ldots,\ell_{n}\left(  \sigma\right)  \right)  <_{\operatorname*{lex}}\left(
\ell_{1}\left(  \tau\right)  ,\ell_{2}\left(  \tau\right)  ,\ldots,\ell
_{n}\left(  \tau\right)  \right)
\]
(according to Definition \ref{def.perm.lehmer.lex-ord}). Hence, Proposition
\ref{prop.perm.lehmer.lex} is proven.
\end{proof}
\end{vershort}

\begin{verlong}
\begin{proof}
[Proof of Proposition \ref{prop.perm.lehmer.lex}.]We have $\left(
\sigma\left(  1\right)  ,\sigma\left(  2\right)  ,\ldots,\sigma\left(
n\right)  \right)  <_{\operatorname*{lex}}\left(  \tau\left(  1\right)
,\tau\left(  2\right)  ,\ldots,\tau\left(  n\right)  \right)  $ if and only if

\begin{itemize}
\item there exists some $k\in\left[  n\right]  $ such that $\sigma\left(
k\right)  \neq\tau\left(  k\right)  $, and the \textbf{smallest} such $k$
satisfies $\sigma\left(  k\right)  <\tau\left(  k\right)  $
\end{itemize}

\noindent(according to Definition \ref{def.perm.lehmer.lex-ord}).

Since we have $\left(  \sigma\left(  1\right)  ,\sigma\left(  2\right)
,\ldots,\sigma\left(  n\right)  \right)  <_{\operatorname*{lex}}\left(
\tau\left(  1\right)  ,\tau\left(  2\right)  ,\ldots,\tau\left(  n\right)
\right)  $ (by assumption), we thus conclude that there exists some
$k\in\left[  n\right]  $ such that $\sigma\left(  k\right)  \neq\tau\left(
k\right)  $, and the \textbf{smallest} such $k$ satisfies $\sigma\left(
k\right)  <\tau\left(  k\right)  $. In other words, the following two
statements hold:

\begin{itemize}
\item There exists some $k\in\left[  n\right]  $ such that $\sigma\left(
k\right)  \neq\tau\left(  k\right)  $.

\item The \textbf{smallest} $k\in\left[  n\right]  $ such that $\sigma\left(
k\right)  \neq\tau\left(  k\right)  $ satisfies $\sigma\left(  k\right)
<\tau\left(  k\right)  $.
\end{itemize}

In particular, there exists some $k\in\left[  n\right]  $ such that
$\sigma\left(  k\right)  \neq\tau\left(  k\right)  $. Let $i$ be the
\textbf{smallest} such $k$. Thus, $i$ is the smallest $k\in\left[  n\right]  $
such that $\sigma\left(  k\right)  \neq\tau\left(  k\right)  $. Hence, $i$ is
a $k\in\left[  n\right]  $ such that $\sigma\left(  i\right)  \neq\tau\left(
i\right)  $. Thus, in particular, $i$ is an element of $\left[  n\right]  $.
Hence, $i\in\left[  n\right]  =\left\{  1,2,\ldots,n\right\}  $ (by the
definition of $\left[  n\right]  $). Hence, $i$ is a positive integer
satisfying $1\leq i\leq n$.

Also, $i$ is the \textbf{smallest} $k\in\left[  n\right]  $ such that
$\sigma\left(  k\right)  \neq\tau\left(  k\right)  $. Thus, every $k\in\left[
n\right]  $ such that $\sigma\left(  k\right)  \neq\tau\left(  k\right)  $
must be $\leq i$. In other words, if $k\in\left[  n\right]  $ is such that
$\sigma\left(  k\right)  \neq\tau\left(  k\right)  $, then%
\begin{equation}
k\geq i. \label{pf.prop.perm.lehmer.lex.long.kgeqi}%
\end{equation}
Thus,
\begin{equation}
\text{each }k\in\left[  i-1\right]  \text{ satisfies }\sigma\left(  k\right)
=\tau\left(  k\right)  \label{pf.prop.perm.lehmer.lex.long.sitau}%
\end{equation}
\footnote{\textit{Proof of (\ref{pf.prop.perm.lehmer.lex.long.sitau}):} Let
$k\in\left[  i-1\right]  $. Thus, $k\in\left[  i-1\right]  =\left\{
1,2,\ldots,i-1\right\}  $ (by the definition of $\left[  i-1\right]  $).
Hence, $1\leq k\leq i-1$. Thus, $k\leq i-1<i\leq n$, so that $k\leq n$ and
therefore $1\leq k\leq n$. Hence, $k\in\left\{  1,2,\ldots,n\right\}  =\left[
n\right]  $.
\par
Assume (for the sake of contradiction) that $\sigma\left(  k\right)  \neq
\tau\left(  k\right)  $. Then, (\ref{pf.prop.perm.lehmer.lex.long.kgeqi})
shows that $k\geq i$. This contradicts $k\leq i-1<i$. This contradiction shows
that our assumption (that $\sigma\left(  k\right)  \neq\tau\left(  k\right)
$) was false. Hence, we don't have $\sigma\left(  k\right)  \neq\tau\left(
k\right)  $. In other words, $\sigma\left(  k\right)  =\tau\left(  k\right)
$. This proves (\ref{pf.prop.perm.lehmer.lex.long.sitau}).}. Thus, Lemma
\ref{lem.perm.lehmer.lex1} \textbf{(b)} shows that we have $\ell_{k}\left(
\sigma\right)  =\ell_{k}\left(  \tau\right)  $ for each $k\in\left[
i-1\right]  $. In other words,%
\begin{equation}
\text{each }k\in\left[  i-1\right]  \text{ satisfies }\ell_{k}\left(
\sigma\right)  =\ell_{k}\left(  \tau\right)  .
\label{pf.prop.perm.lehmer.lex.long.1}%
\end{equation}

We know that the \textbf{smallest} $k\in\left[  n\right]  $ such that
$\sigma\left(  k\right)  \neq\tau\left(  k\right)  $ satisfies $\sigma\left(
k\right)  <\tau\left(  k\right)  $. In other words, $k=i$ satisfies
$\sigma\left(  k\right)  <\tau\left(  k\right)  $ (since $i$ is the
\textbf{smallest} $k\in\left[  n\right]  $ such that $\sigma\left(  k\right)
\neq\tau\left(  k\right)  $). In other words, $\sigma\left(  i\right)
<\tau\left(  i\right)  $. Hence, Lemma \ref{lem.perm.lehmer.lex1} \textbf{(c)}
shows that $\ell_{i}\left(  \sigma\right)  <\ell_{i}\left(  \tau\right)  $.
Thus, $\ell_{i}\left(  \sigma\right)  \neq\ell_{i}\left(  \tau\right)  $. So
we know that $i$ is an element of $\left[  n\right]  $ and satisfies $\ell
_{i}\left(  \sigma\right)  \neq\ell_{i}\left(  \tau\right)  $. In other words,
$i$ is a $k\in\left[  n\right]  $ such that $\ell_{k}\left(  \sigma\right)
\neq\ell_{k}\left(  \tau\right)  $. Hence, there exists some $k\in\left[
n\right]  $ such that $\ell_{k}\left(  \sigma\right)  \neq\ell_{k}\left(
\tau\right)  $ (namely, $k=i$).

Moreover, the \textbf{smallest} $k\in\left[  n\right]  $ such that $\ell
_{k}\left(  \sigma\right)  \neq\ell_{k}\left(  \tau\right)  $ is
$i$\ \ \ \ \footnote{\textit{Proof.} In order to prove this, we need to verify
the following two statements:
\par
\begin{statement}
\textit{Statement 1:} The number $i$ is a $k\in\left[  n\right]  $ such that
$\ell_{k}\left(  \sigma\right)  \neq\ell_{k}\left(  \tau\right)  $.
\end{statement}
\par
\begin{statement}
\textit{Statement 2:} Every $k\in\left[  n\right]  $ such that $\ell
_{k}\left(  \sigma\right)  \neq\ell_{k}\left(  \tau\right)  $ is $\geq i$.
\end{statement}
\par
[\textit{Proof of Statement 1:} We already know that $i$ is a $k\in\left[
n\right]  $ such that $\ell_{k}\left(  \sigma\right)  \neq\ell_{k}\left(
\tau\right)  $. This proves Statement 1.]
\par
[\textit{Proof of Statement 2:} Let $k\in\left[  n\right]  $ be such that
$\ell_{k}\left(  \sigma\right)  \neq\ell_{k}\left(  \tau\right)  $. We shall
prove that $k\geq i$.
\par
Assume the contrary. Thus, $k<i$. But $k\in\left[  n\right]  =\left\{
1,2,\ldots,n\right\}  $. Hence, $k$ is a positive integer. From $k<i$, we
obtain $k\leq i-1$ (since $k$ and $i$ are integers). Thus, $k\in\left\{
1,2,\ldots,i-1\right\}  $ (since $k$ is a positive integer). But the
definition of $\left[  i-1\right]  $ yields $\left[  i-1\right]  =\left\{
1,2,\ldots,i-1\right\}  $. Hence, $k\in\left\{  1,2,\ldots,i-1\right\}
=\left[  i-1\right]  $. Hence, (\ref{pf.prop.perm.lehmer.lex.long.1}) yields
$\ell_{k}\left(  \sigma\right)  =\ell_{k}\left(  \tau\right)  $. This
contradicts $\ell_{k}\left(  \sigma\right)  \neq\ell_{k}\left(  \tau\right)
$.
\par
This contradiction shows that our assumption was wrong. Hence, $k\geq i$ is
proven.
\par
Now, forget that we fixed $k$. We thus have shown that if $k\in\left[
n\right]  $ is such that $\ell_{k}\left(  \sigma\right)  \neq\ell_{k}\left(
\tau\right)  $, then $k\geq i$. In other words, every $k\in\left[  n\right]  $
such that $\ell_{k}\left(  \sigma\right)  \neq\ell_{k}\left(  \tau\right)  $
is $\geq i$. This proves Statement 2.]
\par
Now, both Statement 1 and Statement 2 are proven. Together, these two
statements show that the \textbf{smallest} $k\in\left[  n\right]  $ such that
$\ell_{k}\left(  \sigma\right)  \neq\ell_{k}\left(  \tau\right)  $ is $i$.
Qed.}, and therefore satisfies $\ell_{k}\left(  \sigma\right)  <\ell
_{k}\left(  \tau\right)  $ (because $\ell_{i}\left(  \sigma\right)  <\ell
_{i}\left(  \tau\right)  $).

Hence, we have shown that there exists some $k\in\left[  n\right]  $ such that
$\ell_{k}\left(  \sigma\right)  \neq\ell_{k}\left(  \tau\right)  $, and the
\textbf{smallest} such $k$ satisfies $\ell_{k}\left(  \sigma\right)  <\ell
_{k}\left(  \tau\right)  $.

But we have $\left(  \ell_{1}\left(  \sigma\right)  ,\ell_{2}\left(
\sigma\right)  ,\ldots,\ell_{n}\left(  \sigma\right)  \right)
<_{\operatorname*{lex}}\left(  \ell_{1}\left(  \tau\right)  ,\ell_{2}\left(
\tau\right)  ,\ldots,\ell_{n}\left(  \tau\right)  \right)  $ if and only if

\begin{itemize}
\item there exists some $k\in\left[  n\right]  $ such that $\ell_{k}\left(
\sigma\right)  \neq\ell_{k}\left(  \tau\right)  $, and the \textbf{smallest}
such $k$ satisfies $\ell_{k}\left(  \sigma\right)  <\ell_{k}\left(
\tau\right)  $
\end{itemize}

\noindent(according to Definition \ref{def.perm.lehmer.lex-ord}).

Thus, we have $\left(  \ell_{1}\left(  \sigma\right)  ,\ell_{2}\left(
\sigma\right)  ,\ldots,\ell_{n}\left(  \sigma\right)  \right)
<_{\operatorname*{lex}}\left(  \ell_{1}\left(  \tau\right)  ,\ell_{2}\left(
\tau\right)  ,\ldots,\ell_{n}\left(  \tau\right)  \right)  $ (since there
exists some $k\in\left[  n\right]  $ such that $\ell_{k}\left(  \sigma\right)
\neq\ell_{k}\left(  \tau\right)  $, and the \textbf{smallest} such $k$
satisfies $\ell_{k}\left(  \sigma\right)  <\ell_{k}\left(  \tau\right)  $).
Hence, Proposition \ref{prop.perm.lehmer.lex} is proven.
\end{proof}
\end{verlong}

Next, we prove a further lemma:

\begin{lemma}
\label{lem.perm.lehmer.lex2}Let $\sigma\in S_{n}$ and $\tau\in S_{n}$. Let
$i\in\left[  n\right]  $. Assume that
\[
\text{each }k\in\left[  i-1\right]  \text{ satisfies }\sigma\left(  k\right)
=\tau\left(  k\right)  .
\]
Assume furthermore that $\sigma\left(  i\right)  \neq\tau\left(  i\right)  $.
Then, $\ell_{i}\left(  \sigma\right)  \neq\ell_{i}\left(  \tau\right)  $.
\end{lemma}

\begin{vershort}
\begin{proof}
[Proof of Lemma \ref{lem.perm.lehmer.lex2}.]We have $\sigma\left(  i\right)
\neq\tau\left(  i\right)  $. Hence, we are in one of the following two cases:

\textit{Case 1:} We have $\sigma\left(  i\right)  <\tau\left(  i\right)  $.

\textit{Case 2:} We have $\sigma\left(  i\right)  >\tau\left(  i\right)  $.

Let us first consider Case 1. In this case, we have $\sigma\left(  i\right)
<\tau\left(  i\right)  $. Hence, Lemma \ref{lem.perm.lehmer.lex1} \textbf{(c)}
yields $\ell_{i}\left(  \sigma\right)  <\ell_{i}\left(  \tau\right)  $. Thus,
$\ell_{i}\left(  \sigma\right)  \neq\ell_{i}\left(  \tau\right)  $. This
proves Lemma \ref{lem.perm.lehmer.lex2} in Case 1.

Let us now consider Case 2. In this case, we have $\sigma\left(  i\right)
>\tau\left(  i\right)  $. In other words, $\tau\left(  i\right)
<\sigma\left(  i\right)  $.

But recall that each $k\in\left[  i-1\right]  $ satisfies $\sigma\left(
k\right)  =\tau\left(  k\right)  $. In other words, each $k\in\left[
i-1\right]  $ satisfies $\tau\left(  k\right)  =\sigma\left(  k\right)  $.
Hence, Lemma \ref{lem.perm.lehmer.lex1} \textbf{(c)} (applied to $\tau$ and
$\sigma$ instead of $\sigma$ and $\tau$) yields $\ell_{i}\left(  \tau\right)
<\ell_{i}\left(  \sigma\right)  $. Thus, $\ell_{i}\left(  \sigma\right)
\neq\ell_{i}\left(  \tau\right)  $. Therefore, Lemma
\ref{lem.perm.lehmer.lex2} is proven in Case 2.

We have now proven Lemma \ref{lem.perm.lehmer.lex2} in each of the two Cases 1
and 2. Hence, Lemma \ref{lem.perm.lehmer.lex2} always holds.
\end{proof}
\end{vershort}

\begin{verlong}
\begin{proof}
[Proof of Lemma \ref{lem.perm.lehmer.lex2}.]We have $\sigma\left(  i\right)
\neq\tau\left(  i\right)  $. Hence, we have either $\sigma\left(  i\right)
<\tau\left(  i\right)  $ or $\sigma\left(  i\right)  >\tau\left(  i\right)  $.
Thus, we are in one of the following two cases:

\textit{Case 1:} We have $\sigma\left(  i\right)  <\tau\left(  i\right)  $.

\textit{Case 2:} We have $\sigma\left(  i\right)  >\tau\left(  i\right)  $.

Let us first consider Case 1. In this case, we have $\sigma\left(  i\right)
<\tau\left(  i\right)  $. Hence, Lemma \ref{lem.perm.lehmer.lex1} \textbf{(c)}
yields $\ell_{i}\left(  \sigma\right)  <\ell_{i}\left(  \tau\right)  $. Thus,
$\ell_{i}\left(  \sigma\right)  \neq\ell_{i}\left(  \tau\right)  $. Therefore,
Lemma \ref{lem.perm.lehmer.lex2} is proven in Case 1.

Let us now consider Case 2. In this case, we have $\sigma\left(  i\right)
>\tau\left(  i\right)  $. In other words, $\tau\left(  i\right)
<\sigma\left(  i\right)  $.

But recall that each $k\in\left[  i-1\right]  $ satisfies $\sigma\left(
k\right)  =\tau\left(  k\right)  $. In other words, each $k\in\left[
i-1\right]  $ satisfies $\tau\left(  k\right)  =\sigma\left(  k\right)  $.
Hence, Lemma \ref{lem.perm.lehmer.lex1} \textbf{(c)} (applied to $\tau$ and
$\sigma$ instead of $\sigma$ and $\tau$) yields $\ell_{i}\left(  \tau\right)
<\ell_{i}\left(  \sigma\right)  $. Thus, $\ell_{i}\left(  \tau\right)
\neq\ell_{i}\left(  \sigma\right)  $, so that $\ell_{i}\left(  \sigma\right)
\neq\ell_{i}\left(  \tau\right)  $. Therefore, Lemma
\ref{lem.perm.lehmer.lex2} is proven in Case 2.

We have now proven Lemma \ref{lem.perm.lehmer.lex2} in each of the two Cases 1
and 2. Since these two Cases cover all possibilities, we thus conclude that
Lemma \ref{lem.perm.lehmer.lex2} always holds.
\end{proof}
\end{verlong}

\begin{corollary}
\label{cor.perm.lehmer.inj}Let $\sigma\in S_{n}$ and $\tau\in S_{n}$ such that
$L\left(  \sigma\right)  =L\left(  \tau\right)  $. Then, $\sigma=\tau$.
\end{corollary}

\begin{vershort}
\begin{proof}
[Proof of Corollary \ref{cor.perm.lehmer.inj}.]Assume the contrary (for the
sake of contradiction). Thus, $\sigma\neq\tau$. Hence, there exists some
$k\in\left[  n\right]  $ satisfying $\sigma\left(  k\right)  \neq\tau\left(
k\right)  $.

Let $i$ be the \textbf{smallest} such $k$. Thus, $i$ is an element of $\left[
n\right]  $ and satisfies $\sigma\left(  i\right)  \neq\tau\left(  i\right)
$, but each $k\in\left[  i-1\right]  $ satisfies $\sigma\left(  k\right)
=\tau\left(  k\right)  $. Hence, Lemma \ref{lem.perm.lehmer.lex2} shows that
$\ell_{i}\left(  \sigma\right)  \neq\ell_{i}\left(  \tau\right)  $.

But the definition of $L$ yields $L\left(  \sigma\right)  =\left(  \ell
_{1}\left(  \sigma\right)  ,\ell_{2}\left(  \sigma\right)  ,\ldots,\ell
_{n}\left(  \sigma\right)  \right)  $ and $L\left(  \tau\right)  =\left(
\ell_{1}\left(  \tau\right)  ,\ell_{2}\left(  \tau\right)  ,\ldots,\ell
_{n}\left(  \tau\right)  \right)  $. Thus,%
\[
\left(  \ell_{1}\left(  \sigma\right)  ,\ell_{2}\left(  \sigma\right)
,\ldots,\ell_{n}\left(  \sigma\right)  \right)  =L\left(  \sigma\right)
=L\left(  \tau\right)  =\left(  \ell_{1}\left(  \tau\right)  ,\ell_{2}\left(
\tau\right)  ,\ldots,\ell_{n}\left(  \tau\right)  \right)  .
\]
In other words, each $j\in\left[  n\right]  $ satisfies $\ell_{j}\left(
\sigma\right)  =\ell_{j}\left(  \tau\right)  $. Applying this to $j=i$, we
obtain $\ell_{i}\left(  \sigma\right)  =\ell_{i}\left(  \tau\right)  $. This
contradicts $\ell_{i}\left(  \sigma\right)  \neq\ell_{i}\left(  \tau\right)
$. This contradiction completes our proof of Corollary
\ref{cor.perm.lehmer.inj}.
\end{proof}
\end{vershort}

\begin{verlong}
\begin{proof}
[Proof of Corollary \ref{cor.perm.lehmer.inj}.]We have $\sigma\in S_{n}$. In
other words, $\sigma$ is a permutation of $\left\{  1,2,\ldots,n\right\}  $
(since $S_{n}$ is the set of all permutations of $\left\{  1,2,\ldots
,n\right\}  $). Hence, $\sigma$ is a map $\left\{  1,2,\ldots,n\right\}
\rightarrow\left\{  1,2,\ldots,n\right\}  $. The same argument (applied to
$\tau$ instead of $\sigma$) shows that $\tau$ is a map $\left\{
1,2,\ldots,n\right\}  \rightarrow\left\{  1,2,\ldots,n\right\}  $.

Assume (for the sake of contradiction) that there exists some $k\in\left\{
1,2,\ldots,n\right\}  $ satisfying $\sigma\left(  k\right)  \neq\tau\left(
k\right)  $.

Let $i$ be the \textbf{smallest} such $k$. Thus, $i$ is a $k\in\left\{
1,2,\ldots,n\right\}  $ satisfying $\sigma\left(  k\right)  \neq\tau\left(
k\right)  $. In other words, $i$ is an element of $\left\{  1,2,\ldots
,n\right\}  $ and satisfies $\sigma\left(  i\right)  \neq\tau\left(  i\right)
$.

The definition of $\left[  n\right]  $ yields $\left[  n\right]  =\left\{
1,2,\ldots,n\right\}  $. Thus, $i\in\left\{  1,2,\ldots,n\right\}  =\left[
n\right]  $. Also, from $i\in\left\{  1,2,\ldots,n\right\}  $, we obtain
$i\leq n$.

Also, $i$ is the \textbf{smallest} $k\in\left\{  1,2,\ldots,n\right\}  $
satisfying $\sigma\left(  k\right)  \neq\tau\left(  k\right)  $. Thus, each
$k\in\left\{  1,2,\ldots,n\right\}  $ satisfying $\sigma\left(  k\right)
\neq\tau\left(  k\right)  $ must be $\geq i$. In other words, if $k\in\left\{
1,2,\ldots,n\right\}  $ satisfies $\sigma\left(  k\right)  \neq\tau\left(
k\right)  $, then%
\begin{equation}
k\geq i. \label{pf.cor.perm.lehmer.inj.no-smaller}%
\end{equation}

Thus, each $k\in\left[  i-1\right]  $ satisfies $\sigma\left(  k\right)
=\tau\left(  k\right)  $\ \ \ \ \footnote{\textit{Proof.} Let $k\in\left[
i-1\right]  $. Then, $k\in\left[  i-1\right]  =\left\{  1,2,\ldots
,i-1\right\}  $ (by the definition of $\left[  i-1\right]  $). Hence, $k$ is a
positive integer satisfying $1\leq k\leq i-1$. Now, $k\leq i-1\leq i\leq n$,
so that $k\in\left\{  1,2,\ldots,n\right\}  $ (since $k$ is a positive
integer). Thus, $k\in\left\{  1,2,\ldots,n\right\}  =\left[  n\right]  $.
\par
Assume (for the sake of contradiction) that $\sigma\left(  k\right)  \neq
\tau\left(  k\right)  $. Then, (\ref{pf.cor.perm.lehmer.inj.no-smaller}) shows
that $k\geq i$. This contradicts $k\leq i-1<i$. This contradiction proves that
our assumption (that $\sigma\left(  k\right)  \neq\tau\left(  k\right)  $) was
wrong. Thus, we have $\sigma\left(  k\right)  =\tau\left(  k\right)  $. Qed.}.
Thus, Lemma \ref{lem.perm.lehmer.lex2} shows that $\ell_{i}\left(
\sigma\right)  \neq\ell_{i}\left(  \tau\right)  $.

But the definition of $L$ yields $L\left(  \sigma\right)  =\left(  \ell
_{1}\left(  \sigma\right)  ,\ell_{2}\left(  \sigma\right)  ,\ldots,\ell
_{n}\left(  \sigma\right)  \right)  $ and $L\left(  \tau\right)  =\left(
\ell_{1}\left(  \tau\right)  ,\ell_{2}\left(  \tau\right)  ,\ldots,\ell
_{n}\left(  \tau\right)  \right)  $. Thus,%
\[
\left(  \ell_{1}\left(  \sigma\right)  ,\ell_{2}\left(  \sigma\right)
,\ldots,\ell_{n}\left(  \sigma\right)  \right)  =L\left(  \sigma\right)
=L\left(  \tau\right)  =\left(  \ell_{1}\left(  \tau\right)  ,\ell_{2}\left(
\tau\right)  ,\ldots,\ell_{n}\left(  \tau\right)  \right)  .
\]
In other words, each $j\in\left\{  1,2,\ldots,n\right\}  $ satisfies $\ell
_{j}\left(  \sigma\right)  =\ell_{j}\left(  \tau\right)  $. Applying this to
$j=i$, we obtain $\ell_{i}\left(  \sigma\right)  =\ell_{i}\left(  \tau\right)
$. This contradicts $\ell_{i}\left(  \sigma\right)  \neq\ell_{i}\left(
\tau\right)  $.

This contradiction shows that our assumption (that there exists some
$k\in\left\{  1,2,\ldots,n\right\}  $ satisfying $\sigma\left(  k\right)
\neq\tau\left(  k\right)  $) was false. Hence, there exists no $k\in\left\{
1,2,\ldots,n\right\}  $ satisfying $\sigma\left(  k\right)  \neq\tau\left(
k\right)  $. In other words, each $k\in\left\{  1,2,\ldots,n\right\}  $
satisfies $\sigma\left(  k\right)  =\tau\left(  k\right)  $. In other words,
$\sigma=\tau$ (since $\sigma$ and $\tau$ are maps $\left\{  1,2,\ldots
,n\right\}  \rightarrow\left\{  1,2,\ldots,n\right\}  $). This proves
Corollary \ref{cor.perm.lehmer.inj}.
\end{proof}
\end{verlong}

The next lemma is particularly obvious:

\begin{lemma}
\label{lem.perm.lehmer.sizeH}We have $\left\vert H\right\vert =n!$.
\end{lemma}

\begin{proof}
[Proof of Lemma \ref{lem.perm.lehmer.sizeH}.]For each $k\in\mathbb{N}$, we
have%
\begin{equation}
\left\vert \left[  k\right]  _{0}\right\vert =k+1.
\label{pf.lem.perm.lehmer.sizeH.1}%
\end{equation}

\begin{vershort}
(This follows immediately from $\left[  k\right]  _{0}=\left\{  0,1,\ldots
,k\right\}  $.)
\end{vershort}

\begin{verlong}
[\textit{Proof of (\ref{pf.lem.perm.lehmer.sizeH.1}):} Let $k\in\mathbb{N}$.
The definition of $\left[  k\right]  _{0}$ yields $\left[  k\right]
_{0}=\left\{  0,1,\ldots,k\right\}  $. Hence, $\left\vert \left[  k\right]
_{0}\right\vert =\left\vert \left\{  0,1,\ldots,k\right\}  \right\vert =k+1$
(since $k\in\mathbb{N}$). This proves (\ref{pf.lem.perm.lehmer.sizeH.1}).]
\end{verlong}

But recall that $H=\left[  n-1\right]  _{0}\times\left[  n-2\right]
_{0}\times\cdots\times\left[  n-n\right]  _{0}$. Hence,%
\begin{align*}
\left\vert H\right\vert  &  =\left\vert \left[  n-1\right]  _{0}\times\left[
n-2\right]  _{0}\times\cdots\times\left[  n-n\right]  _{0}\right\vert
=\left\vert \left[  n-1\right]  _{0}\right\vert \cdot\left\vert \left[
n-2\right]  _{0}\right\vert \cdot\cdots\cdot\left\vert \left[  n-n\right]
_{0}\right\vert \\
&  =\prod_{i=1}^{n}\left\vert \left[  n-i\right]  _{0}\right\vert =\prod
_{k=0}^{n-1}\underbrace{\left\vert \left[  k\right]  _{0}\right\vert
}_{\substack{=k+1\\\text{(by (\ref{pf.lem.perm.lehmer.sizeH.1}))}%
}}\ \ \ \ \ \ \ \ \ \ \left(
\begin{array}
[c]{c}%
\text{here, we have substituted }k\text{ for }n-i\\
\text{in the product}%
\end{array}
\right) \\
&  =\prod_{k=0}^{n-1}\left(  k+1\right)  =\prod_{i=1}^{n}%
i\ \ \ \ \ \ \ \ \ \ \left(  \text{here, we have substituted }i\text{ for
}k+1\text{ in the product}\right) \\
&  =1\cdot2\cdot\cdots\cdot n=n!.
\end{align*}
This proves Lemma \ref{lem.perm.lehmer.sizeH}.
\end{proof}

\begin{proof}
[Proof of Theorem \ref{thm.perm.lehmer.bij}.]If $\sigma\in S_{n}$ and $\tau\in
S_{n}$ are such that $L\left(  \sigma\right)  =L\left(  \tau\right)  $, then
$\sigma=\tau$ (by Corollary \ref{cor.perm.lehmer.inj}). In other words, the
map $L$ is injective.

Lemma \ref{lem.perm.lehmer.sizeH} shows that $\left\vert H\right\vert =n!$.
But Corollary \ref{cor.transpos.code.n!} shows that $\left\vert S_{n}%
\right\vert =n!$. Thus, $\left\vert S_{n}\right\vert =n!\geq n!=\left\vert
H\right\vert $. Hence, Lemma \ref{lem.jectivity.pigeon-inj} (applied to
$U=S_{n}$, $V=H$ and $f=L$) shows that we have the following logical
equivalence:%
\[
\left(  L\text{ is injective}\right)  \ \Longleftrightarrow\ \left(  L\text{
is bijective}\right)  .
\]
Hence, $L$ is bijective (since $L$ is injective). In other words, the map
$L:S_{n}\rightarrow H$ is a bijection. This proves Theorem
\ref{thm.perm.lehmer.bij}.
\end{proof}

Finally, it remains to prove Corollary \ref{cor.perm.lehmer.lensum}. We will
need the following fact about products of polynomials:

\begin{lemma}
\label{lem.perm.lehmer.prodrule}For every $i\in\left\{  1,2,\ldots,n\right\}
$, let $Z_{i}$ be a finite set. For every $i\in\left\{  1,2,\ldots,n\right\}
$ and every $k\in Z_{i}$, let $p_{i,k}$ be a polynomial in $x$ with rational
coefficients. Then,%
\[
\prod_{i=1}^{n}\sum_{k\in Z_{i}}p_{i,k}=\sum_{\left(  k_{1},k_{2},\ldots
,k_{n}\right)  \in Z_{1}\times Z_{2}\times\cdots\times Z_{n}}\prod_{i=1}%
^{n}p_{i,k_{i}}.
\]

\end{lemma}

We defer the proof of Lemma \ref{lem.perm.lehmer.prodrule} to a later section:
namely, Lemma \ref{lem.perm.lehmer.prodrule} is the particular case of Lemma
\ref{lem.prodrule.S} (which is proven below) obtained when $\mathbb{K}$ is the
ring $\mathbb{Q}\left[  x\right]  $ (the ring of all polynomials in $x$ with
rational coefficients).

\begin{vershort}
\begin{proof}
[Proof of Corollary \ref{cor.perm.lehmer.lensum}.]Each $m\in\mathbb{N}$
satisfies $\left[  m\right]  _{0}=\left\{  0,1,\ldots,m\right\}  $ (by the
definition of $\left[  m\right]  _{0}$) and thus%
\begin{equation}
\sum_{k\in\left[  m\right]  _{0}}x^{k}=\underbrace{\sum_{k\in\left\{
0,1,\ldots,m\right\}  }}_{=\sum_{k=0}^{m}}x^{k}=\sum_{k=0}^{m}x^{k}.
\label{pf.cor.perm.lehmer.lensum.short.sum1}%
\end{equation}

Each $\sigma\in S_{n}$ satisfies%
\[
\ell\left(  \sigma\right)  =\ell_{1}\left(  \sigma\right)  +\ell_{2}\left(
\sigma\right)  +\cdots+\ell_{n}\left(  \sigma\right)
\]
(by Proposition \ref{prop.perm.lehmer.l}) and therefore%
\begin{equation}
x^{\ell\left(  \sigma\right)  }=x^{\ell_{1}\left(  \sigma\right)  +\ell
_{2}\left(  \sigma\right)  +\cdots+\ell_{n}\left(  \sigma\right)  }%
=x^{\ell_{1}\left(  \sigma\right)  }x^{\ell_{2}\left(  \sigma\right)  }\cdots
x^{\ell_{n}\left(  \sigma\right)  }=\prod_{i=1}^{n}x^{\ell_{i}\left(
\sigma\right)  }. \label{pf.cor.perm.lehmer.lensum.short.monomial}%
\end{equation}

Theorem \ref{thm.perm.lehmer.bij} shows that the map $L$ is a bijection. But
the map $L:S_{n}\rightarrow H$ is defined by
\[
\left(  L\left(  \sigma\right)  =\left(  \ell_{1}\left(  \sigma\right)
,\ell_{2}\left(  \sigma\right)  ,\ldots,\ell_{n}\left(  \sigma\right)
\right)  \ \ \ \ \ \ \ \ \ \ \text{for each }\sigma\in S_{n}\right)  .
\]
Thus, $L$ is the map%
\[
S_{n}\rightarrow H,\ \ \ \ \ \ \ \ \ \ \sigma\mapsto\left(  \ell_{1}\left(
\sigma\right)  ,\ell_{2}\left(  \sigma\right)  ,\ldots,\ell_{n}\left(
\sigma\right)  \right)  .
\]
Hence, the map%
\[
S_{n}\rightarrow H,\ \ \ \ \ \ \ \ \ \ \sigma\mapsto\left(  \ell_{1}\left(
\sigma\right)  ,\ell_{2}\left(  \sigma\right)  ,\ldots,\ell_{n}\left(
\sigma\right)  \right)
\]
is a bijection (since $L$ is a bijection).

We have%
\begin{align}
\sum_{w\in S_{n}}x^{\ell\left(  w\right)  }  &  =\sum_{\sigma\in S_{n}%
}\underbrace{x^{\ell\left(  \sigma\right)  }}_{\substack{=\prod_{i=1}%
^{n}x^{\ell_{i}\left(  \sigma\right)  }\\\text{(by
(\ref{pf.cor.perm.lehmer.lensum.short.monomial}))}}%
}\ \ \ \ \ \ \ \ \ \ \left(  \text{here, we have renamed the summation index
}w\text{ as }\sigma\right) \nonumber\\
&  =\sum_{\sigma\in S_{n}}\prod_{i=1}^{n}x^{\ell_{i}\left(  \sigma\right)
}=\sum_{\left(  k_{1},k_{2},\ldots,k_{n}\right)  \in H}\prod_{i=1}^{n}%
x^{k_{i}}\nonumber\\
&  \ \ \ \ \ \ \ \ \ \ \left(
\begin{array}
[c]{c}%
\text{here, we have substituted }\left(  k_{1},k_{2},\ldots,k_{n}\right) \\
\text{for }\left(  \ell_{1}\left(  \sigma\right)  ,\ell_{2}\left(
\sigma\right)  ,\ldots,\ell_{n}\left(  \sigma\right)  \right)  \text{ in the
sum (since the map}\\
S_{n}\rightarrow H,\ \ \ \ \ \ \ \ \ \ \sigma\mapsto\left(  \ell_{1}\left(
\sigma\right)  ,\ell_{2}\left(  \sigma\right)  ,\ldots,\ell_{n}\left(
\sigma\right)  \right) \\
\text{is a bijection)}%
\end{array}
\right) \nonumber\\
&  =\sum_{\left(  k_{1},k_{2},\ldots,k_{n}\right)  \in\left[  n-1\right]
_{0}\times\left[  n-2\right]  _{0}\times\cdots\times\left[  n-n\right]  _{0}%
}\prod_{i=1}^{n}x^{k_{i}} \label{pf.cor.perm.lehmer.lensum.short.1}%
\end{align}
(since $H=\left[  n-1\right]  _{0}\times\left[  n-2\right]  _{0}\times
\cdots\times\left[  n-n\right]  _{0}$).

But Lemma \ref{lem.perm.lehmer.prodrule} (applied to $Z_{i}=\left[
n-i\right]  _{0}$ and $p_{i,k}=x^{k}$) yields%
\[
\prod_{i=1}^{n}\sum_{k\in\left[  n-i\right]  _{0}}x^{k}=\sum_{\left(
k_{1},k_{2},\ldots,k_{n}\right)  \in\left[  n-1\right]  _{0}\times\left[
n-2\right]  _{0}\times\cdots\times\left[  n-n\right]  _{0}}\prod_{i=1}%
^{n}x^{k_{i}}.
\]
Comparing this with (\ref{pf.cor.perm.lehmer.lensum.short.1}), we obtain%
\begin{align*}
\sum_{w\in S_{n}}x^{\ell\left(  w\right)  }  &  =\prod_{i=1}^{n}\sum
_{k\in\left[  n-i\right]  _{0}}x^{k}=\prod_{m=0}^{n-1}\underbrace{\sum
_{k\in\left[  m\right]  _{0}}x^{k}}_{\substack{=\sum_{k=0}^{m}x^{k}\\\text{(by
(\ref{pf.cor.perm.lehmer.lensum.short.sum1}))}}}\ \ \ \ \ \ \ \ \ \ \left(
\begin{array}
[c]{c}%
\text{here, we have substituted }m\\
\text{for }n-i\text{ in the product}%
\end{array}
\right) \\
&  =\prod_{m=0}^{n-1}\underbrace{\sum_{k=0}^{m}x^{k}}_{=1+x+x^{2}+\cdots
+x^{m}}=\prod_{m=0}^{n-1}\left(  1+x+x^{2}+\cdots+x^{m}\right) \\
&  =1\left(  1+x\right)  \left(  1+x+x^{2}\right)  \cdots\left(
1+x+x^{2}+\cdots+x^{n-1}\right) \\
&  =\left(  1+x\right)  \left(  1+x+x^{2}\right)  \cdots\left(  1+x+x^{2}%
+\cdots+x^{n-1}\right)  .
\end{align*}
This proves Corollary \ref{cor.perm.lehmer.lensum}.
\end{proof}
\end{vershort}

\begin{verlong}
\begin{proof}
[Proof of Corollary \ref{cor.perm.lehmer.lensum}.]Each $m\in\mathbb{N}$
satisfies%
\begin{align}
\sum_{k\in\left[  m\right]  _{0}}x^{k}  &  =\underbrace{\sum_{k\in\left\{
0,1,\ldots,m\right\}  }}_{=\sum_{k=0}^{m}}x^{k}\ \ \ \ \ \ \ \ \ \ \left(
\text{since }\left[  m\right]  _{0}=\left\{  0,1,\ldots,m\right\}  \text{ (by
the definition of }\left[  m\right]  _{0}\text{)}\right) \nonumber\\
&  =\sum_{k=0}^{m}x^{k}. \label{pf.cor.perm.lehmer.lensum.sum1}%
\end{align}

Let $\sigma\in S_{n}$. Then,%
\[
\ell\left(  \sigma\right)  =\ell_{1}\left(  \sigma\right)  +\ell_{2}\left(
\sigma\right)  +\cdots+\ell_{n}\left(  \sigma\right)
\]
(by Proposition \ref{prop.perm.lehmer.l}) and therefore%
\begin{equation}
x^{\ell\left(  \sigma\right)  }=x^{\ell_{1}\left(  \sigma\right)  +\ell
_{2}\left(  \sigma\right)  +\cdots+\ell_{n}\left(  \sigma\right)  }%
=x^{\ell_{1}\left(  \sigma\right)  }x^{\ell_{2}\left(  \sigma\right)  }\cdots
x^{\ell_{n}\left(  \sigma\right)  }=\prod_{i=1}^{n}x^{\ell_{i}\left(
\sigma\right)  }. \label{pf.cor.perm.lehmer.lensum.monomial}%
\end{equation}

Now, forget that we fixed $\sigma$. We thus have proven
(\ref{pf.cor.perm.lehmer.lensum.monomial}) for each $\sigma\in S_{n}$.

Theorem \ref{thm.perm.lehmer.bij} shows that the map $L$ is a bijection.

The map $L:S_{n}\rightarrow H$ is defined by
\[
\left(  L\left(  \sigma\right)  =\left(  \ell_{1}\left(  \sigma\right)
,\ell_{2}\left(  \sigma\right)  ,\ldots,\ell_{n}\left(  \sigma\right)
\right)  \ \ \ \ \ \ \ \ \ \ \text{for each }\sigma\in S_{n}\right)  .
\]
Thus, $L$ is the map%
\[
S_{n}\rightarrow H,\ \ \ \ \ \ \ \ \ \ \sigma\mapsto\left(  \ell_{1}\left(
\sigma\right)  ,\ell_{2}\left(  \sigma\right)  ,\ldots,\ell_{n}\left(
\sigma\right)  \right)  .
\]
Hence, the map%
\[
S_{n}\rightarrow H,\ \ \ \ \ \ \ \ \ \ \sigma\mapsto\left(  \ell_{1}\left(
\sigma\right)  ,\ell_{2}\left(  \sigma\right)  ,\ldots,\ell_{n}\left(
\sigma\right)  \right)
\]
is a bijection (since $L$ is a bijection).

We have%
\begin{align}
\sum_{w\in S_{n}}x^{\ell\left(  w\right)  }  &  =\sum_{\sigma\in S_{n}%
}\underbrace{x^{\ell\left(  \sigma\right)  }}_{\substack{=\prod_{i=1}%
^{n}x^{\ell_{i}\left(  \sigma\right)  }\\\text{(by
(\ref{pf.cor.perm.lehmer.lensum.monomial}))}}}\ \ \ \ \ \ \ \ \ \ \left(
\text{here, we have renamed the summation index }w\text{ as }\sigma\right)
\nonumber\\
&  =\sum_{\sigma\in S_{n}}\prod_{i=1}^{n}x^{\ell_{i}\left(  \sigma\right)
}=\sum_{\left(  k_{1},k_{2},\ldots,k_{n}\right)  \in H}\prod_{i=1}^{n}%
x^{k_{i}}\nonumber\\
&  \ \ \ \ \ \ \ \ \ \ \left(
\begin{array}
[c]{c}%
\text{here, we have substituted }\left(  k_{1},k_{2},\ldots,k_{n}\right) \\
\text{for }\left(  \ell_{1}\left(  \sigma\right)  ,\ell_{2}\left(
\sigma\right)  ,\ldots,\ell_{n}\left(  \sigma\right)  \right) \\
\text{in the sum (since the map}\\
S_{n}\rightarrow H,\ \ \ \ \ \ \ \ \ \ \sigma\mapsto\left(  \ell_{1}\left(
\sigma\right)  ,\ell_{2}\left(  \sigma\right)  ,\ldots,\ell_{n}\left(
\sigma\right)  \right) \\
\text{is a bijection)}%
\end{array}
\right) \nonumber\\
&  =\sum_{\left(  k_{1},k_{2},\ldots,k_{n}\right)  \in\left[  n-1\right]
_{0}\times\left[  n-2\right]  _{0}\times\cdots\times\left[  n-n\right]  _{0}%
}\prod_{i=1}^{n}x^{k_{i}} \label{pf.cor.perm.lehmer.lensum.1}%
\end{align}
(since $H=\left[  n-1\right]  _{0}\times\left[  n-2\right]  _{0}\times
\cdots\times\left[  n-n\right]  _{0}$).

But Lemma \ref{lem.perm.lehmer.prodrule} (applied to $Z_{i}=\left[
n-i\right]  _{0}$ and $p_{i,k}=x^{k}$) yields%
\[
\prod_{i=1}^{n}\sum_{k\in\left[  n-i\right]  _{0}}x^{k}=\sum_{\left(
k_{1},k_{2},\ldots,k_{n}\right)  \in\left[  n-1\right]  _{0}\times\left[
n-2\right]  _{0}\times\cdots\times\left[  n-n\right]  _{0}}\prod_{i=1}%
^{n}x^{k_{i}}.
\]
Comparing this with (\ref{pf.cor.perm.lehmer.lensum.1}), we obtain%
\begin{align*}
\sum_{w\in S_{n}}x^{\ell\left(  w\right)  }  &  =\prod_{i=1}^{n}\sum
_{k\in\left[  n-i\right]  _{0}}x^{k}=\prod_{m=0}^{n-1}\underbrace{\sum
_{k\in\left[  m\right]  _{0}}x^{k}}_{\substack{=\sum_{k=0}^{m}x^{k}\\\text{(by
(\ref{pf.cor.perm.lehmer.lensum.sum1}))}}}\ \ \ \ \ \ \ \ \ \ \left(
\begin{array}
[c]{c}%
\text{here, we have substituted }m\\
\text{for }n-i\text{ in the product}%
\end{array}
\right) \\
&  =\prod_{m=0}^{n-1}\underbrace{\sum_{k=0}^{m}x^{k}}_{\substack{=x^{0}%
+x^{1}+x^{2}+\cdots+x^{m}\\=1+x+x^{2}+\cdots+x^{m}}}=\prod_{m=0}^{n-1}\left(
1+x+x^{2}+\cdots+x^{m}\right) \\
&  =1\cdot\left(  1+x\right)  \cdot\left(  1+x+x^{2}\right)  \cdot\cdots
\cdot\left(  1+x+x^{2}+\cdots+x^{n-1}\right) \\
&  =\left(  1+x\right)  \cdot\left(  1+x+x^{2}\right)  \cdot\cdots\cdot\left(
1+x+x^{2}+\cdots+x^{n-1}\right) \\
&  =\left(  1+x\right)  \left(  1+x+x^{2}\right)  \cdots\left(  1+x+x^{2}%
+\cdots+x^{n-1}\right)  .
\end{align*}
This proves Corollary \ref{cor.perm.lehmer.lensum}.
\end{proof}
\end{verlong}

\subsection{\label{sect.sol.perm.lisitau}Solution to Exercise
\ref{exe.perm.lisitau}}

Throughout Section \ref{sect.sol.perm.lisitau}, we shall use the same
notations that were in use throughout Section \ref{sect.perm.lehmer}. We shall
furthermore use the notation from Definition \ref{def.iverson}. We need
several lemmas to prepare for the solution of Exercise \ref{exe.perm.lisitau}:

\begin{lemma}
\label{lem.sol.perm.lisitau.flipsign}Let $u$ and $v$ be two integers
satisfying $u\neq v$. Then,%
\begin{equation}
\left[  u>v\right]  =1-\left[  u<v\right]
\label{eq.lem.sol.perm.lisitau.flipsign.1}%
\end{equation}
and%
\begin{equation}
\left[  v>u\right]  =\left[  u<v\right]  .
\label{eq.lem.sol.perm.lisitau.flipsign.2}%
\end{equation}

\end{lemma}

\begin{proof}
[Proof of Lemma \ref{lem.sol.perm.lisitau.flipsign}.]The statement $u=v$ is
false (since $u\neq v$). We have the following chain of logical equivalences:%
\begin{align*}
\left(  \text{not }u<v\right)  \  &  \Longleftrightarrow\ \left(  u\geq
v\right)  \ \Longleftrightarrow\ \left(  u>v\text{ or }u=v\right)
\ \Longleftrightarrow\ \left(  u>v\right) \\
&  \ \ \ \ \ \ \ \ \ \ \left(  \text{since the statement }u=v\text{ is
false}\right)  .
\end{align*}
Thus, $\left(  \text{not }u<v\right)  $ and $\left(  u>v\right)  $ are two
equivalent logical statements. Hence, Exercise \ref{exe.iverson-prop}
\textbf{(a)} (applied to $\mathcal{A}=\left(  \text{not }u<v\right)  $ and
$\mathcal{B}=\left(  u>v\right)  $) yields $\left[  \text{not }u<v\right]
=\left[  u>v\right]  $. Hence,%
\[
\left[  u>v\right]  =\left[  \text{not }u<v\right]  =1-\left[  u<v\right]
\]
(by Exercise \ref{exe.iverson-prop} \textbf{(b)} (applied to $\mathcal{A}%
=\left(  u<v\right)  $)). This proves
(\ref{eq.lem.sol.perm.lisitau.flipsign.1}). Thus, it remains to prove
(\ref{eq.lem.sol.perm.lisitau.flipsign.2}).

Clearly, $\left(  v>u\right)  $ and $\left(  u<v\right)  $ are two equivalent
logical statements. Hence, Exercise \ref{exe.iverson-prop} \textbf{(a)}
(applied to $\mathcal{A}=\left(  v>u\right)  $ and $\mathcal{B}=\left(
u<v\right)  $) yields $\left[  v>u\right]  =\left[  u<v\right]  $. Thus,
(\ref{eq.lem.sol.perm.lisitau.flipsign.2}) is proven. This completes the proof
of Lemma \ref{lem.sol.perm.lisitau.flipsign}.
\end{proof}

\begin{lemma}
\label{lem.sol.perm.lisitau.1}Let $n\in\mathbb{N}$. Let $\sigma\in S_{n}$. Let
$i\in\left[  n\right]  $.

\textbf{(a)} We have $\sum_{j\in\left[  n\right]  }\left[  i<j\text{ and
}\sigma\left(  i\right)  >\sigma\left(  j\right)  \right]  =\ell_{i}\left(
\sigma\right)  $.

\textbf{(b)} We have $\sum_{j\in\left[  n\right]  }\left[  i<j\right]  \left[
\sigma\left(  i\right)  >\sigma\left(  j\right)  \right]  =\ell_{i}\left(
\sigma\right)  $.

\textbf{(c)} We have $\sum_{j\in\left[  n\right]  }\left[  i<j\right]  \left(
1-\left[  \sigma\left(  i\right)  <\sigma\left(  j\right)  \right]  \right)
=\ell_{i}\left(  \sigma\right)  $.
\end{lemma}

\begin{vershort}
\begin{proof}
[Proof of Lemma \ref{lem.sol.perm.lisitau.1}.]\textbf{(a)} We have%
\begin{align*}
&  \sum_{j\in\left[  n\right]  }\left[  i<j\text{ and }\sigma\left(  i\right)
>\sigma\left(  j\right)  \right] \\
&  =\sum_{\substack{j\in\left[  n\right]  ;\\i<j\text{ and }\sigma\left(
i\right)  >\sigma\left(  j\right)  }}\underbrace{\left[  i<j\text{ and }%
\sigma\left(  i\right)  >\sigma\left(  j\right)  \right]  }%
_{\substack{=1\\\text{(since }i<j\text{ and }\sigma\left(  i\right)
>\sigma\left(  j\right)  \text{)}}}\\
&  \ \ \ \ \ \ \ \ \ \ +\sum_{\substack{j\in\left[  n\right]  ;\\\text{not
}\left(  i<j\text{ and }\sigma\left(  i\right)  >\sigma\left(  j\right)
\right)  }}\underbrace{\left[  i<j\text{ and }\sigma\left(  i\right)
>\sigma\left(  j\right)  \right]  }_{\substack{=0\\\text{(since we don't have
}\left(  i<j\text{ and }\sigma\left(  i\right)  >\sigma\left(  j\right)
\right)  \text{)}}}\\
&  =\sum_{\substack{j\in\left[  n\right]  ;\\i<j\text{ and }\sigma\left(
i\right)  >\sigma\left(  j\right)  }}1+\underbrace{\sum_{\substack{j\in\left[
n\right]  ;\\\text{not }\left(  i<j\text{ and }\sigma\left(  i\right)
>\sigma\left(  j\right)  \right)  }}0}_{=0}=\sum_{\substack{j\in\left[
n\right]  ;\\i<j\text{ and }\sigma\left(  i\right)  >\sigma\left(  j\right)
}}1\\
&  =\left\vert \left\{  j\in\left[  n\right]  \ \mid\ i<j\text{ and }%
\sigma\left(  i\right)  >\sigma\left(  j\right)  \right\}  \right\vert
\cdot1\\
&  =\left\vert \left\{  j\in\left[  n\right]  \ \mid\ i<j\text{ and }%
\sigma\left(  i\right)  >\sigma\left(  j\right)  \right\}  \right\vert \\
&  =\left(  \text{the number of all }j\in\left[  n\right]  \text{ such that
}i<j\text{ and }\sigma\left(  i\right)  >\sigma\left(  j\right)  \right) \\
&  =\left(  \text{the number of all }j\in\left[  n\right]  \text{ satisfying
}i<j\text{ such that }\sigma\left(  i\right)  >\sigma\left(  j\right)  \right)
\\
&  =\left(  \text{the number of all }j\in\left\{  i+1,i+2,\ldots,n\right\}
\text{ such that }\sigma\left(  i\right)  >\sigma\left(  j\right)  \right) \\
&  \ \ \ \ \ \ \ \ \ \ \left(  \text{since the }j\in\left[  n\right]  \text{
satisfying }i<j\text{ are precisely the }j\in\left\{  i+1,i+2,\ldots
,n\right\}  \right) \\
&  =\ell_{i}\left(  \sigma\right)
\end{align*}
(since $\ell_{i}\left(  \sigma\right)  $ was defined as the number of all
$j\in\left\{  i+1,i+2,\ldots,n\right\}  $ such that $\sigma\left(  i\right)
>\sigma\left(  j\right)  $). This proves Lemma \ref{lem.sol.perm.lisitau.1}
\textbf{(a)}.

\textbf{(b)} Lemma \ref{lem.sol.perm.lisitau.1} \textbf{(a)} yields%
\[
\ell_{i}\left(  \sigma\right)  =\sum_{j\in\left[  n\right]  }%
\underbrace{\left[  i<j\text{ and }\sigma\left(  i\right)  >\sigma\left(
j\right)  \right]  }_{\substack{=\left[  \left(  i<j\right)  \wedge\left(
\sigma\left(  i\right)  >\sigma\left(  j\right)  \right)  \right]  \\=\left[
i<j\right]  \left[  \sigma\left(  i\right)  >\sigma\left(  j\right)  \right]
\\\text{(by Exercise \ref{exe.iverson-prop} \textbf{(c)}}\\\text{(applied to
}\mathcal{A}=\left(  i<j\right)  \text{ and }\mathcal{B}=\left(  \sigma\left(
i\right)  >\sigma\left(  j\right)  \right)  \text{))}}}=\sum_{j\in\left[
n\right]  }\left[  i<j\right]  \left[  \sigma\left(  i\right)  >\sigma\left(
j\right)  \right]  .
\]
This proves Lemma \ref{lem.sol.perm.lisitau.1} \textbf{(b)}.

\textbf{(c)} We shall prove that each $j\in\left[  n\right]  $ satisfies%
\begin{equation}
\left[  i<j\right]  \left[  \sigma\left(  i\right)  >\sigma\left(  j\right)
\right]  =\left[  i<j\right]  \left(  1-\left[  \sigma\left(  i\right)
<\sigma\left(  j\right)  \right]  \right)  .
\label{pf.lem.sol.perm.lisitau.1.short.c.c}%
\end{equation}

[\textit{Proof of (\ref{pf.lem.sol.perm.lisitau.1.short.c.c}):} Let
$j\in\left[  n\right]  $. If we don't have $i<j$, then we have $\left[
i<j\right]  =0$. Hence, if we don't have $i<j$, then
(\ref{pf.lem.sol.perm.lisitau.1.short.c.c}) rewrites as $0=0$, which is
clearly true. Hence, for the rest of this proof of
(\ref{pf.lem.sol.perm.lisitau.1.short.c.c}), we can WLOG assume that we do
have $i<j$. Assume this.

Hence, $i\neq j$. But $\sigma\in S_{n}$. In other words, $\sigma$ is a
permutation of $\left\{  1,2,\ldots,n\right\}  $ (since $S_{n}$ is the set of
all permutations of $\left\{  1,2,\ldots,n\right\}  $). In other words,
$\sigma$ is a bijection $\left\{  1,2,\ldots,n\right\}  \rightarrow\left\{
1,2,\ldots,n\right\}  $. Hence, the map $\sigma$ is bijective, therefore
injective. Hence, from $i\neq j$, we obtain $\sigma\left(  i\right)
\neq\sigma\left(  j\right)  $. Thus, (\ref{eq.lem.sol.perm.lisitau.flipsign.1}%
) (applied to $u=\sigma\left(  i\right)  $ and $v=\sigma\left(  j\right)  $)
yields $\left[  \sigma\left(  i\right)  >\sigma\left(  j\right)  \right]
=1-\left[  \sigma\left(  i\right)  <\sigma\left(  j\right)  \right]  $.
Multiplying both sides of this equality with $\left[  i<j\right]  $, we obtain
$\left[  i<j\right]  \left[  \sigma\left(  i\right)  >\sigma\left(  j\right)
\right]  =\left[  i<j\right]  \left(  1-\left[  \sigma\left(  i\right)
<\sigma\left(  j\right)  \right]  \right)  $. This proves
(\ref{pf.lem.sol.perm.lisitau.1.short.c.c}).]

Now, Lemma \ref{lem.sol.perm.lisitau.1} \textbf{(b)} yields%
\[
\ell_{i}\left(  \sigma\right)  =\sum_{j\in\left[  n\right]  }%
\underbrace{\left[  i<j\right]  \left[  \sigma\left(  i\right)  >\sigma\left(
j\right)  \right]  }_{\substack{=\left[  i<j\right]  \left(  1-\left[
\sigma\left(  i\right)  <\sigma\left(  j\right)  \right]  \right)  \\\text{(by
(\ref{pf.lem.sol.perm.lisitau.1.short.c.c}))}}}=\sum_{j\in\left[  n\right]
}\left[  i<j\right]  \left(  1-\left[  \sigma\left(  i\right)  <\sigma\left(
j\right)  \right]  \right)  .
\]
This proves Lemma \ref{lem.sol.perm.lisitau.1} \textbf{(c)}.
\end{proof}
\end{vershort}

\begin{verlong}
\begin{proof}
[Proof of Lemma \ref{lem.sol.perm.lisitau.1}.]\textbf{(a)} The claim of Lemma
\ref{lem.sol.perm.lisitau.1} \textbf{(a)} is precisely the statement of
(\ref{pf.prop.perm.lehmer.l.li=}), which has been proven during the proof of
Proposition \ref{prop.perm.lehmer.l} (in Section
\ref{sect.sol.perm.lehmer.prove}). Thus, Lemma \ref{lem.sol.perm.lisitau.1}
\textbf{(a)} is already proven.

\textbf{(b)} Lemma \ref{lem.sol.perm.lisitau.1} \textbf{(a)} yields
$\sum_{j\in\left[  n\right]  }\left[  i<j\text{ and }\sigma\left(  i\right)
>\sigma\left(  j\right)  \right]  =\ell_{i}\left(  \sigma\right)  $. Hence,%
\begin{align*}
\ell_{i}\left(  \sigma\right)   &  =\sum_{j\in\left[  n\right]  }%
\underbrace{\left[  i<j\text{ and }\sigma\left(  i\right)  >\sigma\left(
j\right)  \right]  }_{\substack{=\left[  \left(  i<j\right)  \wedge\left(
\sigma\left(  i\right)  >\sigma\left(  j\right)  \right)  \right]
\\\text{(since \textquotedblleft}\wedge\text{\textquotedblright\ is just a
symbol for \textquotedblleft and\textquotedblright)}}}=\sum_{j\in\left[
n\right]  }\underbrace{\left[  \left(  i<j\right)  \wedge\left(  \sigma\left(
i\right)  >\sigma\left(  j\right)  \right)  \right]  }_{\substack{=\left[
i<j\right]  \left[  \sigma\left(  i\right)  >\sigma\left(  j\right)  \right]
\\\text{(by Exercise \ref{exe.iverson-prop} \textbf{(c)}}\\\text{(applied to
}\mathcal{A}=\left(  i<j\right)  \text{ and }\mathcal{B}=\left(  \sigma\left(
i\right)  >\sigma\left(  j\right)  \right)  \text{))}}}\\
&  =\sum_{j\in\left[  n\right]  }\left[  i<j\right]  \left[  \sigma\left(
i\right)  >\sigma\left(  j\right)  \right]  .
\end{align*}
This proves Lemma \ref{lem.sol.perm.lisitau.1} \textbf{(b)}.

\textbf{(c)} We shall prove that each $j\in\left[  n\right]  $ satisfies%
\begin{equation}
\left[  i<j\right]  \left[  \sigma\left(  i\right)  >\sigma\left(  j\right)
\right]  =\left[  i<j\right]  \left(  1-\left[  \sigma\left(  i\right)
<\sigma\left(  j\right)  \right]  \right)  .
\label{pf.lem.sol.perm.lisitau.1.c.c}%
\end{equation}

[\textit{Proof of (\ref{pf.lem.sol.perm.lisitau.1.c.c}):} Let $j\in\left[
n\right]  $. If we don't have $i<j$, then (\ref{pf.lem.sol.perm.lisitau.1.c.c}%
) holds\footnote{\textit{Proof.} Assume that we don't have $i<j$. We must
prove that (\ref{pf.lem.sol.perm.lisitau.1.c.c}) holds.
\par
We don't have $i<j$. Hence, $\left[  i<j\right]  =0$. Thus,
$\underbrace{\left[  i<j\right]  }_{=0}\left[  \sigma\left(  i\right)
>\sigma\left(  j\right)  \right]  =0$. Comparing this with
$\underbrace{\left[  i<j\right]  }_{=0}\left(  1-\left[  \sigma\left(
i\right)  <\sigma\left(  j\right)  \right]  \right)  =0$, we obtain $\left[
i<j\right]  \left[  \sigma\left(  i\right)  >\sigma\left(  j\right)  \right]
=\left[  i<j\right]  \left(  1-\left[  \sigma\left(  i\right)  <\sigma\left(
j\right)  \right]  \right)  $. Hence, (\ref{pf.lem.sol.perm.lisitau.1.c.c})
holds. This completes our proof.}. Hence, for the rest of this proof of
(\ref{pf.lem.sol.perm.lisitau.1.c.c}), we can WLOG assume that we do have
$i<j$. Assume this.

We have $i<j$, thus $i\neq j$. But $\sigma\in S_{n}$. In other words, $\sigma$
is a permutation of $\left\{  1,2,\ldots,n\right\}  $ (since $S_{n}$ is the
set of all permutations of $\left\{  1,2,\ldots,n\right\}  $). In other words,
$\sigma$ is a bijection $\left\{  1,2,\ldots,n\right\}  \rightarrow\left\{
1,2,\ldots,n\right\}  $. Hence, the map $\sigma$ is bijective, therefore
injective. Hence, if we had $\sigma\left(  i\right)  =\sigma\left(  j\right)
$, then we would have $i=j$, which would contradict $i\neq j$. Thus, the
statement $\sigma\left(  i\right)  =\sigma\left(  j\right)  $ is false. In
other words, we have $\sigma\left(  i\right)  \neq\sigma\left(  j\right)  $.
Thus, (\ref{eq.lem.sol.perm.lisitau.flipsign.1}) (applied to $u=\sigma\left(
i\right)  $ and $v=\sigma\left(  j\right)  $) yields $\left[  \sigma\left(
i\right)  >\sigma\left(  j\right)  \right]  =1-\left[  \sigma\left(  i\right)
<\sigma\left(  j\right)  \right]  $. Multiplying both sides of this equality
with $\left[  i<j\right]  $, we obtain $\left[  i<j\right]  \left[
\sigma\left(  i\right)  >\sigma\left(  j\right)  \right]  =\left[  i<j\right]
\left(  1-\left[  \sigma\left(  i\right)  <\sigma\left(  j\right)  \right]
\right)  $. This proves (\ref{pf.lem.sol.perm.lisitau.1.c.c}).]

Now, Lemma \ref{lem.sol.perm.lisitau.1} \textbf{(b)} yields $\sum_{j\in\left[
n\right]  }\left[  i<j\right]  \left[  \sigma\left(  i\right)  >\sigma\left(
j\right)  \right]  =\ell_{i}\left(  \sigma\right)  $. Hence,%
\[
\ell_{i}\left(  \sigma\right)  =\sum_{j\in\left[  n\right]  }%
\underbrace{\left[  i<j\right]  \left[  \sigma\left(  i\right)  >\sigma\left(
j\right)  \right]  }_{\substack{=\left[  i<j\right]  \left(  1-\left[
\sigma\left(  i\right)  <\sigma\left(  j\right)  \right]  \right)  \\\text{(by
(\ref{pf.lem.sol.perm.lisitau.1.c.c}))}}}=\sum_{j\in\left[  n\right]  }\left[
i<j\right]  \left(  1-\left[  \sigma\left(  i\right)  <\sigma\left(  j\right)
\right]  \right)  .
\]
This proves Lemma \ref{lem.sol.perm.lisitau.1} \textbf{(c)}.
\end{proof}
\end{verlong}

\begin{lemma}
\label{lem.sol.perm.lisitau.2}Let $n\in\mathbb{N}$. Let $\sigma\in S_{n}$ and
$\tau\in S_{n}$. Let $i\in\left[  n\right]  $ and $j\in\left[  n\right]  $.
Then,%
\begin{align*}
&  \left[  \tau\left(  i\right)  <\tau\left(  j\right)  \right]  \left(
1-\left[  \sigma\left(  \tau\left(  i\right)  \right)  <\sigma\left(
\tau\left(  j\right)  \right)  \right]  \right)  +\left[  i<j\right]  \left(
1-\left[  \tau\left(  i\right)  <\tau\left(  j\right)  \right]  \right) \\
&  \ \ \ \ \ \ \ \ \ \ -\left[  i<j\right]  \left(  1-\left[  \sigma\left(
\tau\left(  i\right)  \right)  <\sigma\left(  \tau\left(  j\right)  \right)
\right]  \right) \\
&  =\left[  j>i\right]  \left[  \tau\left(  i\right)  >\tau\left(  j\right)
\right]  \left[  \sigma\left(  \tau\left(  j\right)  \right)  >\sigma\left(
\tau\left(  i\right)  \right)  \right] \\
&  \ \ \ \ \ \ \ \ \ \ +\left[  i>j\right]  \left[  \tau\left(  j\right)
>\tau\left(  i\right)  \right]  \left[  \sigma\left(  \tau\left(  i\right)
\right)  >\sigma\left(  \tau\left(  j\right)  \right)  \right]  .
\end{align*}

\end{lemma}

\begin{vershort}
\begin{proof}
[Proof of Lemma \ref{lem.sol.perm.lisitau.2}.]We are in one of the following
two cases:

\textit{Case 1:} We have $i\neq j$.

\textit{Case 2:} We have $i=j$.

Let us first consider Case 1. In this case, we have $i\neq j$.

Both $\tau$ and $\sigma$ belong to $S_{n}$ and thus are permutations of
$\left\{  1,2,\ldots,n\right\}  $ (since $S_{n}$ is the set of all
permutations of $\left\{  1,2,\ldots,n\right\}  $). Hence, both $\tau$ and
$\sigma$ are bijective maps, and thus in particular are injective. From $i\neq
j$, we obtain $\tau\left(  i\right)  \neq\tau\left(  j\right)  $ (since $\tau$
is injective) and therefore $\sigma\left(  \tau\left(  i\right)  \right)
\neq\sigma\left(  \tau\left(  j\right)  \right)  $ (since $\sigma$ is injective).

Define three integers $a$, $b$ and $c$ by%
\[
a=\left[  i<j\right]  ,\ \ \ \ \ \ \ \ \ \ b=\left[  \tau\left(  i\right)
<\tau\left(  j\right)  \right]  \ \ \ \ \ \ \ \ \ \ \text{and}%
\ \ \ \ \ \ \ \ \ \ c=\left[  \sigma\left(  \tau\left(  i\right)  \right)
<\sigma\left(  \tau\left(  j\right)  \right)  \right]  .
\]

Comparing%
\begin{align*}
&  \underbrace{\left[  \tau\left(  i\right)  <\tau\left(  j\right)  \right]
}_{=b}\left(  1-\underbrace{\left[  \sigma\left(  \tau\left(  i\right)
\right)  <\sigma\left(  \tau\left(  j\right)  \right)  \right]  }_{=c}\right)
+\underbrace{\left[  i<j\right]  }_{=a}\left(  1-\underbrace{\left[
\tau\left(  i\right)  <\tau\left(  j\right)  \right]  }_{=b}\right) \\
&  \ \ \ \ \ \ \ \ \ \ -\underbrace{\left[  i<j\right]  }_{=a}\left(
1-\underbrace{\left[  \sigma\left(  \tau\left(  i\right)  \right)
<\sigma\left(  \tau\left(  j\right)  \right)  \right]  }_{=c}\right) \\
&  =b\left(  1-c\right)  +a\left(  1-b\right)  -a\left(  1-c\right)
=b-bc-ab+ac
\end{align*}
with%
\begin{align*}
&  \underbrace{\left[  j>i\right]  }_{\substack{=\left[  i<j\right]
\\\text{(by (\ref{eq.lem.sol.perm.lisitau.flipsign.2})}\\\text{(applied to
}u=i\text{ and }v=j\text{))}}}\ \ \underbrace{\left[  \tau\left(  i\right)
>\tau\left(  j\right)  \right]  }_{\substack{=1-\left[  \tau\left(  i\right)
<\tau\left(  j\right)  \right]  \\\text{(by
(\ref{eq.lem.sol.perm.lisitau.flipsign.1})}\\\text{(applied to }u=\tau\left(
i\right)  \text{ and }v=\tau\left(  j\right)  \text{))}}%
}\ \ \underbrace{\left[  \sigma\left(  \tau\left(  j\right)  \right)
>\sigma\left(  \tau\left(  i\right)  \right)  \right]  }_{\substack{=\left[
\sigma\left(  \tau\left(  i\right)  \right)  <\sigma\left(  \tau\left(
j\right)  \right)  \right]  \\\text{(by
(\ref{eq.lem.sol.perm.lisitau.flipsign.2})}\\\text{(applied to }%
u=\sigma\left(  \tau\left(  i\right)  \right)  \text{ and }v=\sigma\left(
\tau\left(  j\right)  \right)  \text{))}}}\\
&  \ \ \ \ \ \ \ \ \ \ +\underbrace{\left[  i>j\right]  }%
_{\substack{=1-\left[  i<j\right]  \\\text{(by
(\ref{eq.lem.sol.perm.lisitau.flipsign.1})}\\\text{(applied to }u=i\text{ and
}v=j\text{))}}}\ \ \underbrace{\left[  \tau\left(  j\right)  >\tau\left(
i\right)  \right]  }_{\substack{=\left[  \tau\left(  i\right)  <\tau\left(
j\right)  \right]  \\\text{(by (\ref{eq.lem.sol.perm.lisitau.flipsign.2}%
)}\\\text{(applied to }u=\tau\left(  i\right)  \text{ and }v=\tau\left(
j\right)  \text{))}}}\ \ \underbrace{\left[  \sigma\left(  \tau\left(
i\right)  \right)  >\sigma\left(  \tau\left(  j\right)  \right)  \right]
}_{\substack{=1-\left[  \sigma\left(  \tau\left(  i\right)  \right)
<\sigma\left(  \tau\left(  j\right)  \right)  \right]  \\\text{(by
(\ref{eq.lem.sol.perm.lisitau.flipsign.1})}\\\text{(applied to }%
u=\sigma\left(  \tau\left(  i\right)  \right)  \text{ and }v=\sigma\left(
\tau\left(  j\right)  \right)  \text{))}}}\\
&  =\underbrace{\left[  i<j\right]  }_{=a}\left(  1-\underbrace{\left[
\tau\left(  i\right)  <\tau\left(  j\right)  \right]  }_{=b}\right)
\underbrace{\left[  \sigma\left(  \tau\left(  i\right)  \right)
<\sigma\left(  \tau\left(  j\right)  \right)  \right]  }_{=c}\\
&  \ \ \ \ \ \ \ \ \ \ +\left(  1-\underbrace{\left[  i<j\right]  }%
_{=a}\right)  \underbrace{\left[  \tau\left(  i\right)  <\tau\left(  j\right)
\right]  }_{=b}\left(  1-\underbrace{\left[  \sigma\left(  \tau\left(
i\right)  \right)  <\sigma\left(  \tau\left(  j\right)  \right)  \right]
}_{=c}\right) \\
&  =a\left(  1-b\right)  c+\left(  1-a\right)  b\left(  1-c\right)
=b-bc-ab+ac,
\end{align*}
we obtain%
\begin{align*}
&  \left[  \tau\left(  i\right)  <\tau\left(  j\right)  \right]  \left(
1-\left[  \sigma\left(  \tau\left(  i\right)  \right)  <\sigma\left(
\tau\left(  j\right)  \right)  \right]  \right)  +\left[  i<j\right]  \left(
1-\left[  \tau\left(  i\right)  <\tau\left(  j\right)  \right]  \right) \\
&  \ \ \ \ \ \ \ \ \ \ -\left[  i<j\right]  \left(  1-\left[  \sigma\left(
\tau\left(  i\right)  \right)  <\sigma\left(  \tau\left(  j\right)  \right)
\right]  \right) \\
&  =\left[  j>i\right]  \left[  \tau\left(  i\right)  >\tau\left(  j\right)
\right]  \left[  \sigma\left(  \tau\left(  j\right)  \right)  >\sigma\left(
\tau\left(  i\right)  \right)  \right] \\
&  \ \ \ \ \ \ \ \ \ \ +\left[  i>j\right]  \left[  \tau\left(  j\right)
>\tau\left(  i\right)  \right]  \left[  \sigma\left(  \tau\left(  i\right)
\right)  >\sigma\left(  \tau\left(  j\right)  \right)  \right]  .
\end{align*}
Thus, Lemma \ref{lem.sol.perm.lisitau.2} is proven in Case 1.

Let us next consider Case 2. In this case, we have $i=j$. Thus, $i<j$ does not
hold. Hence, $\left[  i<j\right]  =0$. Also, $j>i$ does not hold (since
$j=i$). Hence, $\left[  j>i\right]  =0$. Also, $i>j$ does not hold (since
$i=j$). Hence, $\left[  i>j\right]  =0$. Moreover, $\tau\left(  i\right)
<\tau\left(  j\right)  $ does not hold (since $\tau\left(  \underbrace{i}%
_{=j}\right)  =\tau\left(  j\right)  $). Hence, $\left[  \tau\left(  i\right)
<\tau\left(  j\right)  \right]  =0$. Comparing%
\begin{align*}
&  \underbrace{\left[  \tau\left(  i\right)  <\tau\left(  j\right)  \right]
}_{=0}\left(  1-\left[  \sigma\left(  \tau\left(  i\right)  \right)
<\sigma\left(  \tau\left(  j\right)  \right)  \right]  \right)
+\underbrace{\left[  i<j\right]  }_{=0}\left(  1-\left[  \tau\left(  i\right)
<\tau\left(  j\right)  \right]  \right) \\
&  \ \ \ \ \ \ \ \ \ \ -\underbrace{\left[  i<j\right]  }_{=0}\left(
1-\left[  \sigma\left(  \tau\left(  i\right)  \right)  <\sigma\left(
\tau\left(  j\right)  \right)  \right]  \right) \\
&  =0
\end{align*}
with%
\begin{align*}
&  \underbrace{\left[  j>i\right]  }_{=0}\left[  \tau\left(  i\right)
>\tau\left(  j\right)  \right]  \left[  \sigma\left(  \tau\left(  j\right)
\right)  >\sigma\left(  \tau\left(  i\right)  \right)  \right] \\
&  \ \ \ \ \ \ \ \ \ \ +\underbrace{\left[  i>j\right]  }_{=0}\left[
\tau\left(  j\right)  >\tau\left(  i\right)  \right]  \left[  \sigma\left(
\tau\left(  i\right)  \right)  >\sigma\left(  \tau\left(  j\right)  \right)
\right] \\
&  =0,
\end{align*}
we obtain%
\begin{align*}
&  \left[  \tau\left(  i\right)  <\tau\left(  j\right)  \right]  \left(
1-\left[  \sigma\left(  \tau\left(  i\right)  \right)  <\sigma\left(
\tau\left(  j\right)  \right)  \right]  \right)  +\left[  i<j\right]  \left(
1-\left[  \tau\left(  i\right)  <\tau\left(  j\right)  \right]  \right) \\
&  \ \ \ \ \ \ \ \ \ \ -\left[  i<j\right]  \left(  1-\left[  \sigma\left(
\tau\left(  i\right)  \right)  <\sigma\left(  \tau\left(  j\right)  \right)
\right]  \right) \\
&  =\left[  j>i\right]  \left[  \tau\left(  i\right)  >\tau\left(  j\right)
\right]  \left[  \sigma\left(  \tau\left(  j\right)  \right)  >\sigma\left(
\tau\left(  i\right)  \right)  \right] \\
&  \ \ \ \ \ \ \ \ \ \ +\left[  i>j\right]  \left[  \tau\left(  j\right)
>\tau\left(  i\right)  \right]  \left[  \sigma\left(  \tau\left(  i\right)
\right)  >\sigma\left(  \tau\left(  j\right)  \right)  \right]  .
\end{align*}
Thus, Lemma \ref{lem.sol.perm.lisitau.2} is proven in Case 2.

We have now proven Lemma \ref{lem.sol.perm.lisitau.2} in each of the two Cases
1 and 2. Hence, Lemma \ref{lem.sol.perm.lisitau.2} always holds.
\end{proof}
\end{vershort}

\begin{verlong}
\begin{proof}
[Proof of Lemma \ref{lem.sol.perm.lisitau.2}.]We are in one of the following
two cases:

\textit{Case 1:} We have $i\neq j$.

\textit{Case 2:} We have $i=j$.

Let us first consider Case 1. In this case, we have $i\neq j$.

We have $\tau\in S_{n}$. In other words, $\tau$ is a permutation of $\left\{
1,2,\ldots,n\right\}  $ (since $S_{n}$ is the set of all permutations of
$\left\{  1,2,\ldots,n\right\}  $). In other words, $\tau$ is a bijection
$\left\{  1,2,\ldots,n\right\}  \rightarrow\left\{  1,2,\ldots,n\right\}  $.
Hence, the map $\tau$ is bijective, therefore injective. Hence, if we had
$\tau\left(  i\right)  =\tau\left(  j\right)  $, then we would have $i=j$,
which would contradict $i\neq j$. Thus, the statement $\tau\left(  i\right)
=\tau\left(  j\right)  $ is false. In other words, we have $\tau\left(
i\right)  \neq\tau\left(  j\right)  $.

We have $\sigma\in S_{n}$. In other words, $\sigma$ is a permutation of
$\left\{  1,2,\ldots,n\right\}  $ (since $S_{n}$ is the set of all
permutations of $\left\{  1,2,\ldots,n\right\}  $). In other words, $\sigma$
is a bijection $\left\{  1,2,\ldots,n\right\}  \rightarrow\left\{
1,2,\ldots,n\right\}  $. Hence, the map $\sigma$ is bijective, therefore
injective. Hence, if we had $\sigma\left(  \tau\left(  i\right)  \right)
=\sigma\left(  \tau\left(  j\right)  \right)  $, then we would have
$\tau\left(  i\right)  =\tau\left(  j\right)  $, which would contradict
$\tau\left(  i\right)  \neq\tau\left(  j\right)  $. Thus, the statement
$\sigma\left(  \tau\left(  i\right)  \right)  =\sigma\left(  \tau\left(
j\right)  \right)  $ is false. In other words, we have $\sigma\left(
\tau\left(  i\right)  \right)  \neq\sigma\left(  \tau\left(  j\right)
\right)  $.

Define three integers $a$, $b$ and $c$ by%
\[
a=\left[  i<j\right]  ,\ \ \ \ \ \ \ \ \ \ b=\left[  \tau\left(  i\right)
<\tau\left(  j\right)  \right]  \ \ \ \ \ \ \ \ \ \ \text{and}%
\ \ \ \ \ \ \ \ \ \ c=\left[  \sigma\left(  \tau\left(  i\right)  \right)
<\sigma\left(  \tau\left(  j\right)  \right)  \right]  .
\]
Thus, $\left[  i<j\right]  =a$ and $\left[  \tau\left(  i\right)  <\tau\left(
j\right)  \right]  =b$ and $\left[  \sigma\left(  \tau\left(  i\right)
\right)  <\sigma\left(  \tau\left(  j\right)  \right)  \right]  =c$.

We have $i\neq j$. Hence, (\ref{eq.lem.sol.perm.lisitau.flipsign.1}) (applied
to $u=i$ and $v=j$) yields%
\begin{equation}
\left[  i>j\right]  =1-\underbrace{\left[  i<j\right]  }_{=a}=1-a.
\label{pf.lem.sol.perm.lisitau.2.c1.4}%
\end{equation}
Also, (\ref{eq.lem.sol.perm.lisitau.flipsign.2}) (applied to $u=i$ and $v=j$)
yields
\begin{equation}
\left[  j>i\right]  =\left[  i<j\right]  =a.
\label{pf.lem.sol.perm.lisitau.2.c1.1}%
\end{equation}

Furthermore, $\tau\left(  i\right)  \neq\tau\left(  j\right)  $. Hence,
(\ref{eq.lem.sol.perm.lisitau.flipsign.1}) (applied to $u=\tau\left(
i\right)  $ and $v=\tau\left(  j\right)  $) yields%
\begin{equation}
\left[  \tau\left(  i\right)  >\tau\left(  j\right)  \right]
=1-\underbrace{\left[  \tau\left(  i\right)  <\tau\left(  j\right)  \right]
}_{=b}=1-b. \label{pf.lem.sol.perm.lisitau.2.c1.2}%
\end{equation}
Also, (\ref{eq.lem.sol.perm.lisitau.flipsign.2}) (applied to $u=\tau\left(
i\right)  $ and $v=\tau\left(  j\right)  $) yields%
\begin{equation}
\left[  \tau\left(  j\right)  >\tau\left(  i\right)  \right]  =\left[
\tau\left(  i\right)  <\tau\left(  j\right)  \right]  =b.
\label{pf.lem.sol.perm.lisitau.2.c1.5}%
\end{equation}

Furthermore, $\sigma\left(  \tau\left(  i\right)  \right)  \neq\sigma\left(
\tau\left(  j\right)  \right)  $. Hence,
(\ref{eq.lem.sol.perm.lisitau.flipsign.1}) (applied to $u=\sigma\left(
\tau\left(  i\right)  \right)  $ and $v=\sigma\left(  \tau\left(  j\right)
\right)  $) yields%
\begin{equation}
\left[  \sigma\left(  \tau\left(  i\right)  \right)  >\sigma\left(
\tau\left(  j\right)  \right)  \right]  =1-\underbrace{\left[  \sigma\left(
\tau\left(  i\right)  \right)  <\sigma\left(  \tau\left(  j\right)  \right)
\right]  }_{=c}=1-c. \label{pf.lem.sol.perm.lisitau.2.c1.3}%
\end{equation}
Also, (\ref{eq.lem.sol.perm.lisitau.flipsign.2}) (applied to $u=\sigma\left(
\tau\left(  i\right)  \right)  $ and $v=\sigma\left(  \tau\left(  j\right)
\right)  $) yields%
\begin{equation}
\left[  \sigma\left(  \tau\left(  j\right)  \right)  >\sigma\left(
\tau\left(  i\right)  \right)  \right]  =\left[  \sigma\left(  \tau\left(
i\right)  \right)  <\sigma\left(  \tau\left(  j\right)  \right)  \right]  =c.
\label{pf.lem.sol.perm.lisitau.2.c1.6}%
\end{equation}

Comparing%
\begin{align*}
&  \underbrace{\left[  \tau\left(  i\right)  <\tau\left(  j\right)  \right]
}_{=b}\left(  1-\underbrace{\left[  \sigma\left(  \tau\left(  i\right)
\right)  <\sigma\left(  \tau\left(  j\right)  \right)  \right]  }_{=c}\right)
+\underbrace{\left[  i<j\right]  }_{=a}\left(  1-\underbrace{\left[
\tau\left(  i\right)  <\tau\left(  j\right)  \right]  }_{=b}\right) \\
&  \ \ \ \ \ \ \ \ \ \ -\underbrace{\left[  i<j\right]  }_{=a}\left(
1-\underbrace{\left[  \sigma\left(  \tau\left(  i\right)  \right)
<\sigma\left(  \tau\left(  j\right)  \right)  \right]  }_{=c}\right) \\
&  =b\left(  1-c\right)  +a\left(  1-b\right)  -a\left(  1-c\right)
=b-bc+a-ab-\left(  a-ac\right) \\
&  =b-bc-ab+ac
\end{align*}
with%
\begin{align*}
&  \underbrace{\left[  j>i\right]  }_{\substack{=a\\\text{(by
(\ref{pf.lem.sol.perm.lisitau.2.c1.1}))}}}\underbrace{\left[  \tau\left(
i\right)  >\tau\left(  j\right)  \right]  }_{\substack{=1-b\\\text{(by
(\ref{pf.lem.sol.perm.lisitau.2.c1.2}))}}}\underbrace{\left[  \sigma\left(
\tau\left(  j\right)  \right)  >\sigma\left(  \tau\left(  i\right)  \right)
\right]  }_{\substack{=c\\\text{(by (\ref{pf.lem.sol.perm.lisitau.2.c1.6}))}%
}}\\
&  \ \ \ \ \ \ \ \ \ \ +\underbrace{\left[  i>j\right]  }%
_{\substack{=1-a\\\text{(by (\ref{pf.lem.sol.perm.lisitau.2.c1.4}))}%
}}\underbrace{\left[  \tau\left(  j\right)  >\tau\left(  i\right)  \right]
}_{\substack{=b\\\text{(by (\ref{pf.lem.sol.perm.lisitau.2.c1.5}))}%
}}\underbrace{\left[  \sigma\left(  \tau\left(  i\right)  \right)
>\sigma\left(  \tau\left(  j\right)  \right)  \right]  }%
_{\substack{=1-c\\\text{(by (\ref{pf.lem.sol.perm.lisitau.2.c1.3}))}}}\\
&  =\underbrace{a\left(  1-b\right)  c}_{=ac-abc}+\underbrace{\left(
1-a\right)  b\left(  1-c\right)  }_{=b-ab-bc+abc}=\left(  ac-abc\right)
+\left(  b-ab-bc+abc\right) \\
&  =b-bc-ab+ac,
\end{align*}
we obtain%
\begin{align*}
&  \left[  \tau\left(  i\right)  <\tau\left(  j\right)  \right]  \left(
1-\left[  \sigma\left(  \tau\left(  i\right)  \right)  <\sigma\left(
\tau\left(  j\right)  \right)  \right]  \right)  +\left[  i<j\right]  \left(
1-\left[  \tau\left(  i\right)  <\tau\left(  j\right)  \right]  \right) \\
&  \ \ \ \ \ \ \ \ \ \ -\left[  i<j\right]  \left(  1-\left[  \sigma\left(
\tau\left(  i\right)  \right)  <\sigma\left(  \tau\left(  j\right)  \right)
\right]  \right) \\
&  =\left[  j>i\right]  \left[  \tau\left(  i\right)  >\tau\left(  j\right)
\right]  \left[  \sigma\left(  \tau\left(  j\right)  \right)  >\sigma\left(
\tau\left(  i\right)  \right)  \right] \\
&  \ \ \ \ \ \ \ \ \ \ +\left[  i>j\right]  \left[  \tau\left(  j\right)
>\tau\left(  i\right)  \right]  \left[  \sigma\left(  \tau\left(  i\right)
\right)  >\sigma\left(  \tau\left(  j\right)  \right)  \right]  .
\end{align*}
Thus, Lemma \ref{lem.sol.perm.lisitau.2} is proven in Case 1.

Let us next consider Case 2. In this case, we have $i=j$. Thus, $i<j$ does not
hold. Hence, $\left[  i<j\right]  =0$. Also, $j>i$ does not hold (since
$j=i$). Hence, $\left[  j>i\right]  =0$. Also, $i>j$ does not hold (since
$i=j$). Hence, $\left[  i>j\right]  =0$. Moreover, $\tau\left(  i\right)
<\tau\left(  j\right)  $ does not hold (since $\tau\left(  \underbrace{i}%
_{=j}\right)  =\tau\left(  j\right)  $). Hence, $\left[  \tau\left(  i\right)
<\tau\left(  j\right)  \right]  =0$. Comparing%
\begin{align*}
&  \underbrace{\left[  \tau\left(  i\right)  <\tau\left(  j\right)  \right]
}_{=0}\left(  1-\left[  \sigma\left(  \tau\left(  i\right)  \right)
<\sigma\left(  \tau\left(  j\right)  \right)  \right]  \right)
+\underbrace{\left[  i<j\right]  }_{=0}\left(  1-\left[  \tau\left(  i\right)
<\tau\left(  j\right)  \right]  \right) \\
&  \ \ \ \ \ \ \ \ \ \ -\underbrace{\left[  i<j\right]  }_{=0}\left(
1-\left[  \sigma\left(  \tau\left(  i\right)  \right)  <\sigma\left(
\tau\left(  j\right)  \right)  \right]  \right) \\
&  =0\left(  1-\left[  \sigma\left(  \tau\left(  i\right)  \right)
<\sigma\left(  \tau\left(  j\right)  \right)  \right]  \right)  +0\left(
1-\left[  \tau\left(  i\right)  <\tau\left(  j\right)  \right]  \right)
-0\left(  1-\left[  \sigma\left(  \tau\left(  i\right)  \right)
<\sigma\left(  \tau\left(  j\right)  \right)  \right]  \right) \\
&  =0
\end{align*}
with%
\begin{align*}
&  \underbrace{\left[  j>i\right]  }_{=0}\left[  \tau\left(  i\right)
>\tau\left(  j\right)  \right]  \left[  \sigma\left(  \tau\left(  j\right)
\right)  >\sigma\left(  \tau\left(  i\right)  \right)  \right] \\
&  \ \ \ \ \ \ \ \ \ \ +\underbrace{\left[  i>j\right]  }_{=0}\left[
\tau\left(  j\right)  >\tau\left(  i\right)  \right]  \left[  \sigma\left(
\tau\left(  i\right)  \right)  >\sigma\left(  \tau\left(  j\right)  \right)
\right] \\
&  =0\left[  \tau\left(  i\right)  >\tau\left(  j\right)  \right]  \left[
\sigma\left(  \tau\left(  j\right)  \right)  >\sigma\left(  \tau\left(
i\right)  \right)  \right]  +0\left[  \tau\left(  j\right)  >\tau\left(
i\right)  \right]  \left[  \sigma\left(  \tau\left(  i\right)  \right)
>\sigma\left(  \tau\left(  j\right)  \right)  \right] \\
&  =0,
\end{align*}
we obtain%
\begin{align*}
&  \left[  \tau\left(  i\right)  <\tau\left(  j\right)  \right]  \left(
1-\left[  \sigma\left(  \tau\left(  i\right)  \right)  <\sigma\left(
\tau\left(  j\right)  \right)  \right]  \right)  +\left[  i<j\right]  \left(
1-\left[  \tau\left(  i\right)  <\tau\left(  j\right)  \right]  \right) \\
&  \ \ \ \ \ \ \ \ \ \ -\left[  i<j\right]  \left(  1-\left[  \sigma\left(
\tau\left(  i\right)  \right)  <\sigma\left(  \tau\left(  j\right)  \right)
\right]  \right) \\
&  =\left[  j>i\right]  \left[  \tau\left(  i\right)  >\tau\left(  j\right)
\right]  \left[  \sigma\left(  \tau\left(  j\right)  \right)  >\sigma\left(
\tau\left(  i\right)  \right)  \right] \\
&  \ \ \ \ \ \ \ \ \ \ +\left[  i>j\right]  \left[  \tau\left(  j\right)
>\tau\left(  i\right)  \right]  \left[  \sigma\left(  \tau\left(  i\right)
\right)  >\sigma\left(  \tau\left(  j\right)  \right)  \right]  .
\end{align*}
Thus, Lemma \ref{lem.sol.perm.lisitau.2} is proven in Case 2.

We have now proven Lemma \ref{lem.sol.perm.lisitau.2} in each of the two Cases
1 and 2. Since these two Cases cover all possibilities, we thus conclude that
Lemma \ref{lem.sol.perm.lisitau.2} always holds.
\end{proof}
\end{verlong}

Let us now recall Definition \ref{def.perm.si}.

\begin{lemma}
\label{lem.sol.perm.lisitau.sk0}Let $n\in\mathbb{N}$. Let $k\in\left\{
1,2,\ldots,n-1\right\}  $. Let $i\in\left[  n\right]  $. Then, $s_{k}\left(
i\right)  =i+\left[  i=k\right]  -\left[  i=k+1\right]  $.
\end{lemma}

\begin{vershort}
\begin{proof}
[Proof of Lemma \ref{lem.sol.perm.lisitau.sk0}.]Recall that $s_{k}$ is the
permutation in $S_{n}$ that swaps $k$ with $k+1$ but leaves all other numbers
unchanged (by the definition of $s_{k}$). Thus, we have $s_{k}\left(
k\right)  =k+1$ and $s_{k}\left(  k+1\right)  =k$ and
\begin{equation}
\left(  s_{k}\left(  j\right)  =j\ \ \ \ \ \ \ \ \ \ \text{for each }%
j\in\left\{  1,2,\ldots,n\right\}  \setminus\left\{  k,k+1\right\}  \right)  .
\label{pf.lem.sol.perm.lisitau.sk0.short.1}%
\end{equation}

Now, we are in one of the following three cases:

\textit{Case 1:} We have $i=k$.

\textit{Case 2:} We have $i=k+1$.

\textit{Case 3:} We have neither $i=k$ nor $i=k+1$.

Let us first consider Case 1. In this case, we have $i=k$. Thus, $s_{k}\left(
\underbrace{i}_{=k}\right)  =s_{k}\left(  k\right)  =k+1$. Comparing this with%
\begin{align*}
i+\left[  i=k\right]  -\left[  i=k+1\right]   &  =k+\underbrace{\left[
k=k\right]  }_{\substack{=1\\\text{(since }k=k\text{)}}}-\underbrace{\left[
k=k+1\right]  }_{\substack{=0\\\text{(since }k=k+1\text{ is false)}%
}}\ \ \ \ \ \ \ \ \ \ \left(  \text{since }i=k\right) \\
&  =k+1-0=k+1,
\end{align*}
we obtain $s_{k}\left(  i\right)  =i+\left[  i=k\right]  -\left[
i=k+1\right]  $. Hence, Lemma \ref{lem.sol.perm.lisitau.sk0} is proven in Case 1.

An analogous argument proves Lemma \ref{lem.sol.perm.lisitau.sk0} in Case 2.

Let us finally consider Case 3. In this case, we have neither $i=k$ nor
$i=k+1$. In other words, $i\notin\left\{  k,k+1\right\}  $. Combining
$i\in\left[  n\right]  =\left\{  1,2,\ldots,n\right\}  $ with $i\notin\left\{
k,k+1\right\}  $, we obtain $i\in\left\{  1,2,\ldots,n\right\}  \setminus
\left\{  k,k+1\right\}  $. Thus, (\ref{pf.lem.sol.perm.lisitau.sk0.short.1})
(applied to $j=i$) yields $s_{k}\left(  i\right)  =i$. Comparing this with%
\[
i+\underbrace{\left[  i=k\right]  }_{\substack{=0\\\text{(since we don't have
}i=k\text{)}}}-\underbrace{\left[  i=k+1\right]  }_{\substack{=0\\\text{(since
we don't have }i=k+1\text{)}}}=i+0-0=i,
\]
we obtain $s_{k}\left(  i\right)  =i+\left[  i=k\right]  -\left[
i=k+1\right]  $. Hence, Lemma \ref{lem.sol.perm.lisitau.sk0} is proven in Case 3.

We have now proven Lemma \ref{lem.sol.perm.lisitau.sk0} in each of the three
Cases 1, 2 and 3. Hence, Lemma \ref{lem.sol.perm.lisitau.sk0} always holds.
\end{proof}
\end{vershort}

\begin{verlong}
\begin{proof}
[Proof of Lemma \ref{lem.sol.perm.lisitau.sk0}.]Recall that $s_{k}$ is the
permutation in $S_{n}$ that swaps $k$ with $k+1$ but leaves all other numbers
unchanged (by the definition of $s_{k}$). Thus, we have $s_{k}\left(
k\right)  =k+1$ and $s_{k}\left(  k+1\right)  =k$ and
\begin{equation}
\left(  s_{k}\left(  j\right)  =j\ \ \ \ \ \ \ \ \ \ \text{for each }%
j\in\left\{  1,2,\ldots,n\right\}  \setminus\left\{  k,k+1\right\}  \right)  .
\label{pf.lem.sol.perm.lisitau.sk0.1}%
\end{equation}

Now, we are in one of the following three cases:

\textit{Case 1:} We have $i=k$.

\textit{Case 2:} We have $i=k+1$.

\textit{Case 3:} We have neither $i=k$ nor $i=k+1$.

Let us first consider Case 1. In this case, we have $i=k$. Thus, $s_{k}\left(
\underbrace{i}_{=k}\right)  =s_{k}\left(  k\right)  =k+1$. Comparing this with%
\begin{align*}
i+\left[  i=k\right]  -\left[  i=k+1\right]   &  =k+\underbrace{\left[
k=k\right]  }_{\substack{=1\\\text{(since }k=k\text{)}}}-\underbrace{\left[
k=k+1\right]  }_{\substack{=0\\\text{(since }k=k+1\text{ is false}%
\\\text{(since }k\neq k+1\text{))}}}\ \ \ \ \ \ \ \ \ \ \left(  \text{since
}i=k\right) \\
&  =k+1-0=k+1,
\end{align*}
we obtain $s_{k}\left(  i\right)  =i+\left[  i=k\right]  -\left[
i=k+1\right]  $. Hence, Lemma \ref{lem.sol.perm.lisitau.sk0} is proven in Case 1.

Let us next consider Case 2. In this case, we have $i=k+1$. Thus,
$s_{k}\left(  \underbrace{i}_{=k+1}\right)  =s_{k}\left(  k+1\right)  =k$.
Comparing this with%
\begin{align*}
i+\left[  i=k\right]  -\left[  i=k+1\right]   &  =\left(  k+1\right)
+\underbrace{\left[  k+1=k\right]  }_{\substack{=0\\\text{(since }k+1=k\text{
is false}\\\text{(since }k+1\neq k\text{))}}}-\underbrace{\left[
k+1=k+1\right]  }_{\substack{=1\\\text{(since }k+1=k+1\text{)}}%
}\ \ \ \ \ \ \ \ \ \ \left(  \text{since }i=k+1\right) \\
&  =\left(  k+1\right)  +0-1=k,
\end{align*}
we obtain $s_{k}\left(  i\right)  =i+\left[  i=k\right]  -\left[
i=k+1\right]  $. Hence, Lemma \ref{lem.sol.perm.lisitau.sk0} is proven in Case 2.

Let us finally consider Case 3. In this case, we have neither $i=k$ nor
$i=k+1$. Thus, $i\notin\left\{  k,k+1\right\}  $%
\ \ \ \ \footnote{\textit{Proof.} Assume the contrary. Thus, $i\in\left\{
k,k+1\right\}  $. Hence, $i$ equals either $k$ or $k+1$. In other words, we
have either $i=k$ or $i=k+1$. This contradicts the fact that we have neither
$i=k$ nor $i=k+1$. This contradiction shows that our assumption was wrong,
qed.}. But $i\in\left[  n\right]  =\left\{  1,2,\ldots,n\right\}  $ (by the
definition of $\left[  n\right]  $). Combining this with $i\notin\left\{
k,k+1\right\}  $, we obtain $i\in\left\{  1,2,\ldots,n\right\}  \setminus
\left\{  k,k+1\right\}  $. Thus, (\ref{pf.lem.sol.perm.lisitau.sk0.1})
(applied to $j=i$) yields $s_{k}\left(  i\right)  =i$. Comparing this with%
\[
i+\underbrace{\left[  i=k\right]  }_{\substack{=0\\\text{(since we don't have
}i=k\text{)}}}-\underbrace{\left[  i=k+1\right]  }_{\substack{=0\\\text{(since
we don't have }i=k+1\text{)}}}=i+0-0=i,
\]
we obtain $s_{k}\left(  i\right)  =i+\left[  i=k\right]  -\left[
i=k+1\right]  $. Hence, Lemma \ref{lem.sol.perm.lisitau.sk0} is proven in Case 3.

We have now proven Lemma \ref{lem.sol.perm.lisitau.sk0} in each of the three
Cases 1, 2 and 3. Since these three Cases cover all possibilities, we thus
conclude that Lemma \ref{lem.sol.perm.lisitau.sk0} always holds.
\end{proof}
\end{verlong}

\begin{lemma}
\label{lem.sol.perm.lisitau.sk1}Let $n\in\mathbb{N}$. Let $k\in\left\{
1,2,\ldots,n-1\right\}  $. Let $i\in\left[  n\right]  $ and $j\in\left[
n\right]  $. Then:

\textbf{(a)} If $j>i$, then $\left[  s_{k}\left(  i\right)  >s_{k}\left(
j\right)  \right]  =\left[  i=k\right]  \left[  j=k+1\right]  $.

\textbf{(b)} If $j\leq i$, then $\left[  i=k\right]  \left[  j=k+1\right]  =0$.

\textbf{(c)} We have%
\[
\left[  j>i\right]  \left[  s_{k}\left(  i\right)  >s_{k}\left(  j\right)
\right]  =\left[  i=k\right]  \left[  j=k+1\right]  .
\]

\end{lemma}

\begin{vershort}
\begin{proof}
[Proof of Lemma \ref{lem.sol.perm.lisitau.sk1}.]The definition of $s_{k}$
yields $s_{k}\left(  k\right)  =k+1$ and $s_{k}\left(  k+1\right)  =k$.

\textbf{(a)} Assume that $j>i$. Thus, $i<j$, so that $i\leq j-1$ (since both
$i$ and $j$ are integers). In other words, $i+1\leq j$.

Every statement $\mathcal{A}$ satisfies $\left[  \mathcal{A}\right]
\in\left\{  0,1\right\}  $ and thus $\left[  \mathcal{A}\right]  \geq0$ and
$\left[  \mathcal{A}\right]  \leq1$. These facts lead to $\left[
i=k+1\right]  \geq0$, $\left[  j=k\right]  \geq0$, $\left[  j=k+1\right]
\leq1$ and $\left[  i=k\right]  \leq1$.

We are in one of the following three cases:

\textit{Case 1:} We have $i\neq k$.

\textit{Case 2:} We have $j\neq k+1$.

\textit{Case 3:} We have neither $i\neq k$ nor $j\neq k+1$.

(Of course, there can be overlap between Case 1 and Case 2.)

Let us first consider Case 1. In this case, we have $i\neq k$. In other words,
we don't have $i=k$. Hence, $\left[  i=k\right]  =0$. Now, Lemma
\ref{lem.sol.perm.lisitau.sk0} yields%
\[
s_{k}\left(  i\right)  =i+\underbrace{\left[  i=k\right]  }_{=0}%
-\underbrace{\left[  i=k+1\right]  }_{\geq0}\leq i+0-0=i\leq j-1.
\]

Also, Lemma \ref{lem.sol.perm.lisitau.sk0} (applied to $j$ instead of $i$)
yields%
\[
s_{k}\left(  j\right)  =j+\underbrace{\left[  j=k\right]  }_{\geq
0}-\underbrace{\left[  j=k+1\right]  }_{\leq1}\geq j+0-1=j-1.
\]
Thus, $j-1\leq s_{k}\left(  j\right)  $, so that $s_{k}\left(  i\right)  \leq
j-1\leq s_{k}\left(  j\right)  $. Thus, we cannot have $s_{k}\left(  i\right)
>s_{k}\left(  j\right)  $. Hence, $\left[  s_{k}\left(  i\right)
>s_{k}\left(  j\right)  \right]  =0$. Comparing this with $\underbrace{\left[
i=k\right]  }_{=0}\left[  j=k+1\right]  =0$, we obtain $\left[  s_{k}\left(
i\right)  >s_{k}\left(  j\right)  \right]  =\left[  i=k\right]  \left[
j=k+1\right]  $. Hence, Lemma \ref{lem.sol.perm.lisitau.sk1} \textbf{(a)} is
proven in Case 1.

Let us next consider Case 2. In this case, we have $j\neq k+1$. In other
words, we don't have $j=k+1$. Hence, $\left[  j=k+1\right]  =0$. Now, Lemma
\ref{lem.sol.perm.lisitau.sk0} (applied to $j$ instead of $i$) yields%
\[
s_{k}\left(  j\right)  =j+\underbrace{\left[  j=k\right]  }_{\geq
0}-\underbrace{\left[  j=k+1\right]  }_{=0}\geq j+0-0=j.
\]
In other words, $j\leq s_{k}\left(  j\right)  $.

Also, Lemma \ref{lem.sol.perm.lisitau.sk0} yields%
\[
s_{k}\left(  i\right)  =i+\underbrace{\left[  i=k\right]  }_{\leq
1}-\underbrace{\left[  i=k+1\right]  }_{\geq0}\leq i+1-0=i+1\leq j\leq
s_{k}\left(  j\right)  .
\]
Thus, we cannot have $s_{k}\left(  i\right)  >s_{k}\left(  j\right)  $. Hence,
$\left[  s_{k}\left(  i\right)  >s_{k}\left(  j\right)  \right]  =0$.
Comparing this with $\left[  i=k\right]  \underbrace{\left[  j=k+1\right]
}_{=0}=0$, we obtain $\left[  s_{k}\left(  i\right)  >s_{k}\left(  j\right)
\right]  =\left[  i=k\right]  \left[  j=k+1\right]  $. Hence, Lemma
\ref{lem.sol.perm.lisitau.sk1} \textbf{(a)} is proven in Case 2.

Let us finally consider Case 3. In this case, we have neither $i\neq k$ nor
$j\neq k+1$. Thus, we have both $i=k$ and $j=k+1$. Now, $s_{k}\left(
\underbrace{i}_{=k}\right)  =s_{k}\left(  k\right)  =k+1$ and $s_{k}\left(
\underbrace{j}_{=k+1}\right)  =s_{k}\left(  k+1\right)  =k$, so that
$s_{k}\left(  i\right)  =k+1>k=s_{k}\left(  j\right)  $. Hence, $\left[
s_{k}\left(  i\right)  >s_{k}\left(  j\right)  \right]  =1$. Comparing this
with $\underbrace{\left[  i=k\right]  }_{\substack{=1\\\text{(since
}i=k\text{)}}}\underbrace{\left[  j=k+1\right]  }_{\substack{=1\\\text{(since
}j=k+1\text{)}}}=1$, we obtain $\left[  s_{k}\left(  i\right)  >s_{k}\left(
j\right)  \right]  =\left[  i=k\right]  \left[  j=k+1\right]  $. Hence, Lemma
\ref{lem.sol.perm.lisitau.sk1} \textbf{(a)} is proven in Case 3.

We have now proven Lemma \ref{lem.sol.perm.lisitau.sk1} \textbf{(a)} in each
of the three Cases 1, 2 and 3. Hence, Lemma \ref{lem.sol.perm.lisitau.sk1}
\textbf{(a)} always holds.

\textbf{(b)} Assume that $j\leq i$. We are in one of the following two cases:

\textit{Case 1:} We have $i\leq k$.

\textit{Case 2:} We have $i>k$.

Let us first consider Case 1. In this case, we have $i\leq k$. Thus, $j\leq
i\leq k<k+1$. Hence, we don't have $j=k+1$. Thus, $\left[  j=k+1\right]  =0$.
Thus, $\left[  i=k\right]  \underbrace{\left[  j=k+1\right]  }_{=0}=0$. Hence,
Lemma \ref{lem.sol.perm.lisitau.sk1} \textbf{(b)} is proven in Case 1.

Let us now consider Case 2. In this case, we have $i>k$. Hence, we don't have
$i=k$. Thus, $\left[  i=k\right]  =0$. Thus, $\underbrace{\left[  i=k\right]
}_{=0}\left[  j=k+1\right]  =0$. Hence, Lemma \ref{lem.sol.perm.lisitau.sk1}
\textbf{(b)} is proven in Case 2.

We have now proven Lemma \ref{lem.sol.perm.lisitau.sk1} \textbf{(b)} in each
of the two Cases 1 and 2. Thus, Lemma \ref{lem.sol.perm.lisitau.sk1}
\textbf{(b)} always holds.

\textbf{(c)} We are in one of the following two cases:

\textit{Case 1:} We have $j>i$.

\textit{Case 2:} We don't have $j>i$.

Let us first consider Case 1. In this case, we have $j>i$. Hence, $\left[
j>i\right]  =1$. Hence, $\underbrace{\left[  j>i\right]  }_{=1}\left[
s_{k}\left(  i\right)  >s_{k}\left(  j\right)  \right]  =\left[  s_{k}\left(
i\right)  >s_{k}\left(  j\right)  \right]  =\left[  i=k\right]  \left[
j=k+1\right]  $ (by Lemma \ref{lem.sol.perm.lisitau.sk1} \textbf{(a)}). Thus,
Lemma \ref{lem.sol.perm.lisitau.sk1} \textbf{(c)} is proven in Case 1.

Let us next consider Case 2. In this case, we don't have $j>i$. Hence,
$\left[  j>i\right]  =0$. Also, $j\leq i$ (since we don't have $j>i$); thus,
Lemma \ref{lem.sol.perm.lisitau.sk1} \textbf{(b)} yields $\left[  i=k\right]
\left[  j=k+1\right]  =0$. Comparing this with $\underbrace{\left[
j>i\right]  }_{=0}\left[  s_{k}\left(  i\right)  >s_{k}\left(  j\right)
\right]  =0$, we obtain $\left[  j>i\right]  \left[  s_{k}\left(  i\right)
>s_{k}\left(  j\right)  \right]  =\left[  i=k\right]  \left[  j=k+1\right]  $.
Therefore, Lemma \ref{lem.sol.perm.lisitau.sk1} \textbf{(c)} is proven in Case 2.

We have now proven Lemma \ref{lem.sol.perm.lisitau.sk1} \textbf{(c)} in each
of the two Cases 1 and 2. Thus, Lemma \ref{lem.sol.perm.lisitau.sk1}
\textbf{(c)} always holds.
\end{proof}
\end{vershort}

\begin{verlong}
\begin{proof}
[Proof of Lemma \ref{lem.sol.perm.lisitau.sk1}.]Recall that $s_{k}$ is the
permutation in $S_{n}$ that swaps $k$ with $k+1$ but leaves all other numbers
unchanged (by the definition of $s_{k}$). Thus, we have $s_{k}\left(
k\right)  =k+1$ and $s_{k}\left(  k+1\right)  =k$ and
\[
\left(  s_{k}\left(  j\right)  =j\ \ \ \ \ \ \ \ \ \ \text{for each }%
j\in\left\{  1,2,\ldots,n\right\}  \setminus\left\{  k,k+1\right\}  \right)
.
\]

\textbf{(a)} Assume that $j>i$. Thus, $i<j$, so that $i\leq j-1$ (since both
$i$ and $j$ are integers). In other words, $i+1\leq j$.

We are in one of the following three cases:

\textit{Case 1:} We have $i\neq k$.

\textit{Case 2:} We have $j\neq k+1$.

\textit{Case 3:} We have neither $i\neq k$ nor $j\neq k+1$.

(Of course, there can be overlap between Case 1 and Case 2.)

Let us first consider Case 1. In this case, we have $i\neq k$. In other words,
we don't have $i=k$. Hence, $\left[  i=k\right]  =0$. But recall that every
statement $\mathcal{A}$ satisfies $\left[  \mathcal{A}\right]  \in\left\{
0,1\right\}  $ and thus $\left[  \mathcal{A}\right]  \geq0$. Applying this to
$\mathcal{A}=\left(  i=k+1\right)  $, we obtain $\left[  i=k+1\right]  \geq0$.
Now, Lemma \ref{lem.sol.perm.lisitau.sk0} yields%
\[
s_{k}\left(  i\right)  =i+\underbrace{\left[  i=k\right]  }_{=0}%
-\underbrace{\left[  i=k+1\right]  }_{\geq0}\leq i+0-0=i\leq j-1.
\]

Also, recall that every statement $\mathcal{A}$ satisfies $\left[
\mathcal{A}\right]  \in\left\{  0,1\right\}  $ and thus $\left[
\mathcal{A}\right]  \leq1$. Applying this to $\mathcal{A}=\left(
j=k+1\right)  $, we obtain $\left[  j=k+1\right]  \leq1$. Furthermore, every
statement $\mathcal{A}$ satisfies $\left[  \mathcal{A}\right]  \in\left\{
0,1\right\}  $ and thus $\left[  \mathcal{A}\right]  \geq0$. Applying this to
$\mathcal{A}=\left(  j=k\right)  $, we obtain $\left[  j=k\right]  \geq0$.
Now, Lemma \ref{lem.sol.perm.lisitau.sk0} (applied to $j$ instead of $i$)
yields%
\[
s_{k}\left(  j\right)  =j+\underbrace{\left[  j=k\right]  }_{\geq
0}-\underbrace{\left[  j=k+1\right]  }_{\leq1}\geq j+0-1=j-1.
\]
Thus, $j-1\leq s_{k}\left(  j\right)  $, so that $s_{k}\left(  i\right)  \leq
j-1\leq s_{k}\left(  j\right)  $. Thus, we cannot have $s_{k}\left(  i\right)
>s_{k}\left(  j\right)  $. Hence, $\left[  s_{k}\left(  i\right)
>s_{k}\left(  j\right)  \right]  =0$. Comparing this with $\underbrace{\left[
i=k\right]  }_{=0}\left[  j=k+1\right]  =0$, we obtain $\left[  s_{k}\left(
i\right)  >s_{k}\left(  j\right)  \right]  =\left[  i=k\right]  \left[
j=k+1\right]  $. Hence, Lemma \ref{lem.sol.perm.lisitau.sk1} \textbf{(a)} is
proven in Case 1.

Let us next consider Case 2. In this case, we have $j\neq k+1$. In other
words, we don't have $j=k+1$. Hence, $\left[  j=k+1\right]  =0$. But recall
that every statement $\mathcal{A}$ satisfies $\left[  \mathcal{A}\right]
\in\left\{  0,1\right\}  $ and thus $\left[  \mathcal{A}\right]  \geq0$.
Applying this to $\mathcal{A}=\left(  j=k\right)  $, we obtain $\left[
j=k\right]  \geq0$. Now, Lemma \ref{lem.sol.perm.lisitau.sk0} (applied to $j$
instead of $i$) yields%
\[
s_{k}\left(  j\right)  =j+\underbrace{\left[  j=k\right]  }_{\geq
0}-\underbrace{\left[  j=k+1\right]  }_{=0}\geq j+0-0=j.
\]
In other words, $j\leq s_{k}\left(  j\right)  $.

Also, recall that every statement $\mathcal{A}$ satisfies $\left[
\mathcal{A}\right]  \in\left\{  0,1\right\}  $ and thus $\left[
\mathcal{A}\right]  \leq1$. Applying this to $\mathcal{A}=\left(  i=k\right)
$, we obtain $\left[  i=k\right]  \leq1$. Furthermore, every statement
$\mathcal{A}$ satisfies $\left[  \mathcal{A}\right]  \in\left\{  0,1\right\}
$ and thus $\left[  \mathcal{A}\right]  \geq0$. Applying this to
$\mathcal{A}=\left(  i=k+1\right)  $, we obtain $\left[  i=k+1\right]  \geq0$.
Now, Lemma \ref{lem.sol.perm.lisitau.sk0} yields%
\[
s_{k}\left(  i\right)  =i+\underbrace{\left[  i=k\right]  }_{\leq
1}-\underbrace{\left[  i=k+1\right]  }_{\geq0}\leq i+1-0=i+1\leq j\leq
s_{k}\left(  j\right)  .
\]
Thus, we cannot have $s_{k}\left(  i\right)  >s_{k}\left(  j\right)  $. Hence,
$\left[  s_{k}\left(  i\right)  >s_{k}\left(  j\right)  \right]  =0$.
Comparing this with $\left[  i=k\right]  \underbrace{\left[  j=k+1\right]
}_{=0}=0$, we obtain $\left[  s_{k}\left(  i\right)  >s_{k}\left(  j\right)
\right]  =\left[  i=k\right]  \left[  j=k+1\right]  $. Hence, Lemma
\ref{lem.sol.perm.lisitau.sk1} \textbf{(a)} is proven in Case 2.

Let us finally consider Case 3. In this case, we have neither $i\neq k$ nor
$j\neq k+1$. Thus, we have both $i=k$ and $j=k+1$. Now, $s_{k}\left(
\underbrace{i}_{=k}\right)  =s_{k}\left(  k\right)  =k+1$ and $s_{k}\left(
\underbrace{j}_{=k+1}\right)  =s_{k}\left(  k+1\right)  =k$, so that
$s_{k}\left(  i\right)  =k+1>k=s_{k}\left(  j\right)  $. Hence, $\left[
s_{k}\left(  i\right)  >s_{k}\left(  j\right)  \right]  =1$. Comparing this
with $\underbrace{\left[  i=k\right]  }_{\substack{=1\\\text{(since
}i=k\text{)}}}\underbrace{\left[  j=k+1\right]  }_{\substack{=1\\\text{(since
}j=k+1\text{)}}}=1$, we obtain $\left[  s_{k}\left(  i\right)  >s_{k}\left(
j\right)  \right]  =\left[  i=k\right]  \left[  j=k+1\right]  $. Hence, Lemma
\ref{lem.sol.perm.lisitau.sk1} \textbf{(a)} is proven in Case 3.

We have now proven Lemma \ref{lem.sol.perm.lisitau.sk1} \textbf{(a)} in each
of the three Cases 1, 2 and 3. Since these three Cases cover all
possibilities, we thus conclude that Lemma \ref{lem.sol.perm.lisitau.sk1}
\textbf{(a)} always holds.

\textbf{(b)} Assume that $j\leq i$. We are in one of the following two cases:

\textit{Case 1:} We have $i\leq k$.

\textit{Case 2:} We have $i>k$.

Let us first consider Case 1. In this case, we have $i\leq k$. Thus, $j\leq
i\leq k<k+1$, so that $j\neq k+1$. Hence, we don't have $j=k+1$. Thus,
$\left[  j=k+1\right]  =0$. Thus, $\left[  i=k\right]  \underbrace{\left[
j=k+1\right]  }_{=0}=0$. Hence, Lemma \ref{lem.sol.perm.lisitau.sk1}
\textbf{(b)} is proven in Case 1.

Let us now consider Case 2. In this case, we have $i>k$. Thus, $i\neq k$.
Hence, we don't have $i=k$. Thus, $\left[  i=k\right]  =0$. Thus,
$\underbrace{\left[  i=k\right]  }_{=0}\left[  j=k+1\right]  =0$. Hence, Lemma
\ref{lem.sol.perm.lisitau.sk1} \textbf{(b)} is proven in Case 2.

We have now proven Lemma \ref{lem.sol.perm.lisitau.sk1} \textbf{(b)} in each
of the two Cases 1 and 2. Since these two Cases cover all possibilities, we
thus conclude that Lemma \ref{lem.sol.perm.lisitau.sk1} \textbf{(b)} always holds.

\textbf{(c)} We are in one of the following two cases:

\textit{Case 1:} We have $j>i$.

\textit{Case 2:} We don't have $j>i$.

Let us first consider Case 1. In this case, we have $j>i$. Hence, $\left[
j>i\right]  =1$. Hence, $\underbrace{\left[  j>i\right]  }_{=1}\left[
s_{k}\left(  i\right)  >s_{k}\left(  j\right)  \right]  =\left[  s_{k}\left(
i\right)  >s_{k}\left(  j\right)  \right]  =\left[  i=k\right]  \left[
j=k+1\right]  $ (by Lemma \ref{lem.sol.perm.lisitau.sk1} \textbf{(a)}). Thus,
Lemma \ref{lem.sol.perm.lisitau.sk1} \textbf{(c)} is proven in Case 1.

Let us next consider Case 2. In this case, we don't have $j>i$. Hence,
$\left[  j>i\right]  =0$. Also, $j\leq i$ (since we don't have $j>i$); thus,
Lemma \ref{lem.sol.perm.lisitau.sk1} \textbf{(b)} yields $\left[  i=k\right]
\left[  j=k+1\right]  =0$. Comparing this with $\underbrace{\left[
j>i\right]  }_{=0}\left[  s_{k}\left(  i\right)  >s_{k}\left(  j\right)
\right]  =0$, we obtain $\left[  j>i\right]  \left[  s_{k}\left(  i\right)
>s_{k}\left(  j\right)  \right]  =\left[  i=k\right]  \left[  j=k+1\right]  $.
Therefore, Lemma \ref{lem.sol.perm.lisitau.sk1} \textbf{(c)} is proven in Case 2.

We have now proven Lemma \ref{lem.sol.perm.lisitau.sk1} \textbf{(c)} in each
of the two Cases 1 and 2. Since these two Cases cover all possibilities, we
thus conclude that Lemma \ref{lem.sol.perm.lisitau.sk1} \textbf{(c)} always holds.
\end{proof}
\end{verlong}

\begin{lemma}
\label{lem.sol.perm.lisitau.trivsum}Let $P$ be a finite set. Let $p\in P$.
Then, $\sum_{j\in P}\left[  j=p\right]  =1$.
\end{lemma}

\begin{proof}
[Proof of Lemma \ref{lem.sol.perm.lisitau.trivsum}.]We have $p\in P$. Thus, we
can split off the addend for $j=p$ from the sum $\sum_{j\in P}\left[
j=p\right]  $. We thus obtain%
\[
\sum_{j\in P}\left[  j=p\right]  =\underbrace{\left[  p=p\right]
}_{\substack{=1\\\text{(since }p=p\text{)}}}+\sum_{\substack{j\in P;\\j\neq
p}}\underbrace{\left[  j=p\right]  }_{\substack{=0\\\text{(since we don't have
}j=p\\\text{(since }j\neq p\text{))}}}=1+\underbrace{\sum_{\substack{j\in
P;\\j\neq p}}0}_{=0}=1.
\]
This proves Lemma \ref{lem.sol.perm.lisitau.trivsum}.
\end{proof}

\begin{proof}
[Solution to Exercise \ref{exe.perm.lisitau}.]We have $\left[  n\right]
=\left\{  1,2,\ldots,n\right\}  $ (by the definition of $\left[  n\right]  $).

\begin{vershort}
We have $\tau\in S_{n}$. In other words, $\tau$ is a permutation of $\left\{
1,2,\ldots,n\right\}  $. In other words, $\tau$ is a permutation of $\left[
n\right]  $ (since $\left[  n\right]  =\left\{  1,2,\ldots,n\right\}  $). In
other words, $\tau$ is a bijection $\left[  n\right]  \rightarrow\left[
n\right]  $.
\end{vershort}

\begin{verlong}
We have $\tau\in S_{n}$. In other words, $\tau$ is a permutation of $\left\{
1,2,\ldots,n\right\}  $ (since $S_{n}$ is the set of all permutations of
$\left\{  1,2,\ldots,n\right\}  $). In other words, $\tau$ is a bijection
$\left\{  1,2,\ldots,n\right\}  \rightarrow\left\{  1,2,\ldots,n\right\}  $.
In other words, $\tau$ is a bijection $\left[  n\right]  \rightarrow\left[
n\right]  $ (since $\left[  n\right]  =\left\{  1,2,\ldots,n\right\}  $). In
other words, the map $\tau:\left[  n\right]  \rightarrow\left[  n\right]  $ is
a bijection.
\end{verlong}

\textbf{(a)} Let $i\in\left[  n\right]  $. Then, $\tau\left(  i\right)
\in\left[  n\right]  $ (since $\tau$ is a bijection $\left[  n\right]
\rightarrow\left[  n\right]  $). Hence, Lemma \ref{lem.sol.perm.lisitau.1}
\textbf{(c)} (applied to $\tau\left(  i\right)  $ instead of $i$) yields%
\[
\sum_{j\in\left[  n\right]  }\left[  \tau\left(  i\right)  <j\right]  \left(
1-\left[  \sigma\left(  \tau\left(  i\right)  \right)  <\sigma\left(
j\right)  \right]  \right)  =\ell_{\tau\left(  i\right)  }\left(
\sigma\right)  .
\]
Hence,%
\begin{align}
\ell_{\tau\left(  i\right)  }\left(  \sigma\right)   &  =\sum_{j\in\left[
n\right]  }\left[  \tau\left(  i\right)  <j\right]  \left(  1-\left[
\sigma\left(  \tau\left(  i\right)  \right)  <\sigma\left(  j\right)  \right]
\right) \nonumber\\
&  =\sum_{j\in\left[  n\right]  }\left[  \tau\left(  i\right)  <\tau\left(
j\right)  \right]  \left(  1-\left[  \sigma\left(  \tau\left(  i\right)
\right)  <\sigma\left(  \tau\left(  j\right)  \right)  \right]  \right)
\label{sol.perm.lisitau.a.1}%
\end{align}
(here, we have substituted $\tau\left(  j\right)  $ for $j$ in the sum, since
the map $\tau:\left[  n\right]  \rightarrow\left[  n\right]  $ is a bijection).

Also, Lemma \ref{lem.sol.perm.lisitau.1} \textbf{(c)} (applied to $\sigma
\circ\tau$ instead of $\sigma$) yields%
\[
\sum_{j\in\left[  n\right]  }\left[  i<j\right]  \left(  1-\left[  \left(
\sigma\circ\tau\right)  \left(  i\right)  <\left(  \sigma\circ\tau\right)
\left(  j\right)  \right]  \right)  =\ell_{i}\left(  \sigma\circ\tau\right)
.
\]
Hence,%
\begin{align}
\ell_{i}\left(  \sigma\circ\tau\right)   &  =\sum_{j\in\left[  n\right]
}\left[  i<j\right]  \left(  1-\left[  \underbrace{\left(  \sigma\circ
\tau\right)  \left(  i\right)  }_{=\sigma\left(  \tau\left(  i\right)
\right)  }<\underbrace{\left(  \sigma\circ\tau\right)  \left(  j\right)
}_{=\sigma\left(  \tau\left(  j\right)  \right)  }\right]  \right) \nonumber\\
&  =\sum_{j\in\left[  n\right]  }\left[  i<j\right]  \left(  1-\left[
\sigma\left(  \tau\left(  i\right)  \right)  <\sigma\left(  \tau\left(
j\right)  \right)  \right]  \right)  . \label{sol.perm.lisitau.a.2}%
\end{align}

Moreover, Lemma \ref{lem.sol.perm.lisitau.1} \textbf{(c)} (applied to $\tau$
instead of $\sigma$) yields%
\begin{equation}
\sum_{j\in\left[  n\right]  }\left[  i<j\right]  \left(  1-\left[  \tau\left(
i\right)  <\tau\left(  j\right)  \right]  \right)  =\ell_{i}\left(
\tau\right)  . \label{sol.perm.lisitau.a.3}%
\end{equation}

But Lemma \ref{lem.sol.perm.lisitau.2} shows that every $j\in\left[  n\right]
$ satisfies%
\begin{align*}
&  \left[  \tau\left(  i\right)  <\tau\left(  j\right)  \right]  \left(
1-\left[  \sigma\left(  \tau\left(  i\right)  \right)  <\sigma\left(
\tau\left(  j\right)  \right)  \right]  \right)  +\left[  i<j\right]  \left(
1-\left[  \tau\left(  i\right)  <\tau\left(  j\right)  \right]  \right) \\
&  \ \ \ \ \ \ \ \ \ \ -\left[  i<j\right]  \left(  1-\left[  \sigma\left(
\tau\left(  i\right)  \right)  <\sigma\left(  \tau\left(  j\right)  \right)
\right]  \right) \\
&  =\left[  j>i\right]  \left[  \tau\left(  i\right)  >\tau\left(  j\right)
\right]  \left[  \sigma\left(  \tau\left(  j\right)  \right)  >\sigma\left(
\tau\left(  i\right)  \right)  \right] \\
&  \ \ \ \ \ \ \ \ \ \ +\left[  i>j\right]  \left[  \tau\left(  j\right)
>\tau\left(  i\right)  \right]  \left[  \sigma\left(  \tau\left(  i\right)
\right)  >\sigma\left(  \tau\left(  j\right)  \right)  \right]  .
\end{align*}
Summing up these equalities for all $j\in\left[  n\right]  $, we obtain%
\begin{align}
&  \sum_{j\in\left[  n\right]  }\left(  \left[  \tau\left(  i\right)
<\tau\left(  j\right)  \right]  \left(  1-\left[  \sigma\left(  \tau\left(
i\right)  \right)  <\sigma\left(  \tau\left(  j\right)  \right)  \right]
\right)  +\left[  i<j\right]  \left(  1-\left[  \tau\left(  i\right)
<\tau\left(  j\right)  \right]  \right)  \vphantom{\dfrac{1}{1}}\right.
\nonumber\\
&  \ \ \ \ \ \ \ \ \ \ \ \ \ \ \ \ \ \ \ \ \left.
\vphantom{\dfrac{1}{1}}-\left[  i<j\right]  \left(  1-\left[  \sigma\left(
\tau\left(  i\right)  \right)  <\sigma\left(  \tau\left(  j\right)  \right)
\right]  \right)  \right) \nonumber\\
&  =\sum_{j\in\left[  n\right]  }\left(  \left[  j>i\right]  \left[
\tau\left(  i\right)  >\tau\left(  j\right)  \right]  \left[  \sigma\left(
\tau\left(  j\right)  \right)  >\sigma\left(  \tau\left(  i\right)  \right)
\right]  \vphantom{\dfrac{1}{1}}\right. \nonumber\\
&  \ \ \ \ \ \ \ \ \ \ \ \ \ \ \ \ \ \ \ \ \left.
\vphantom{\dfrac{1}{1}}+\left[  i>j\right]  \left[  \tau\left(  j\right)
>\tau\left(  i\right)  \right]  \left[  \sigma\left(  \tau\left(  i\right)
\right)  >\sigma\left(  \tau\left(  j\right)  \right)  \right]  \right)  .
\label{sol.perm.lisitau.a.5}%
\end{align}

Now,%
\begin{align*}
&  \underbrace{\ell_{\tau\left(  i\right)  }\left(  \sigma\right)
}_{\substack{=\sum_{j\in\left[  n\right]  }\left[  \tau\left(  i\right)
<\tau\left(  j\right)  \right]  \left(  1-\left[  \sigma\left(  \tau\left(
i\right)  \right)  <\sigma\left(  \tau\left(  j\right)  \right)  \right]
\right)  \\\text{(by (\ref{sol.perm.lisitau.a.1}))}}}+\underbrace{\ell
_{i}\left(  \tau\right)  }_{\substack{=\sum_{j\in\left[  n\right]  }\left[
i<j\right]  \left(  1-\left[  \tau\left(  i\right)  <\tau\left(  j\right)
\right]  \right)  \\\text{(by (\ref{sol.perm.lisitau.a.3}))}}%
}-\underbrace{\ell_{i}\left(  \sigma\circ\tau\right)  }_{\substack{=\sum
_{j\in\left[  n\right]  }\left[  i<j\right]  \left(  1-\left[  \sigma\left(
\tau\left(  i\right)  \right)  <\sigma\left(  \tau\left(  j\right)  \right)
\right]  \right)  \\\text{(by (\ref{sol.perm.lisitau.a.2}))}}}\\
&  =\sum_{j\in\left[  n\right]  }\left[  \tau\left(  i\right)  <\tau\left(
j\right)  \right]  \left(  1-\left[  \sigma\left(  \tau\left(  i\right)
\right)  <\sigma\left(  \tau\left(  j\right)  \right)  \right]  \right)
+\sum_{j\in\left[  n\right]  }\left[  i<j\right]  \left(  1-\left[
\tau\left(  i\right)  <\tau\left(  j\right)  \right]  \right) \\
&  \ \ \ \ \ \ \ \ \ \ -\sum_{j\in\left[  n\right]  }\left[  i<j\right]
\left(  1-\left[  \sigma\left(  \tau\left(  i\right)  \right)  <\sigma\left(
\tau\left(  j\right)  \right)  \right]  \right) \\
&  =\sum_{j\in\left[  n\right]  }\left(  \left[  \tau\left(  i\right)
<\tau\left(  j\right)  \right]  \left(  1-\left[  \sigma\left(  \tau\left(
i\right)  \right)  <\sigma\left(  \tau\left(  j\right)  \right)  \right]
\right)  +\left[  i<j\right]  \left(  1-\left[  \tau\left(  i\right)
<\tau\left(  j\right)  \right]  \right)  \vphantom{\dfrac{1}{1}}\right. \\
&  \ \ \ \ \ \ \ \ \ \ \ \ \ \ \ \ \ \ \ \ \left.
\vphantom{\dfrac{1}{1}}-\left[  i<j\right]  \left(  1-\left[  \sigma\left(
\tau\left(  i\right)  \right)  <\sigma\left(  \tau\left(  j\right)  \right)
\right]  \right)  \right) \\
&  =\sum_{j\in\left[  n\right]  }\left(  \left[  j>i\right]  \left[
\tau\left(  i\right)  >\tau\left(  j\right)  \right]  \left[  \sigma\left(
\tau\left(  j\right)  \right)  >\sigma\left(  \tau\left(  i\right)  \right)
\right]  \vphantom{\dfrac{1}{1}}\right. \\
&  \ \ \ \ \ \ \ \ \ \ \ \ \ \ \ \ \ \ \ \ \left.
\vphantom{\dfrac{1}{1}}+\left[  i>j\right]  \left[  \tau\left(  j\right)
>\tau\left(  i\right)  \right]  \left[  \sigma\left(  \tau\left(  i\right)
\right)  >\sigma\left(  \tau\left(  j\right)  \right)  \right]  \right) \\
&  \ \ \ \ \ \ \ \ \ \ \left(  \text{by (\ref{sol.perm.lisitau.a.5})}\right)
\\
&  =\sum_{j\in\left[  n\right]  }\left[  j>i\right]  \left[  \tau\left(
i\right)  >\tau\left(  j\right)  \right]  \left[  \sigma\left(  \tau\left(
j\right)  \right)  >\sigma\left(  \tau\left(  i\right)  \right)  \right] \\
&  \ \ \ \ \ \ \ \ \ \ +\sum_{j\in\left[  n\right]  }\left[  i>j\right]
\left[  \tau\left(  j\right)  >\tau\left(  i\right)  \right]  \left[
\sigma\left(  \tau\left(  i\right)  \right)  >\sigma\left(  \tau\left(
j\right)  \right)  \right]  .
\end{align*}
This solves Exercise \ref{exe.perm.lisitau} \textbf{(a)}.

\textbf{(b)} We have $\left[  n\right]  =\left\{  1,2,\ldots,n\right\}  $ and
thus $\sum_{i\in\left[  n\right]  }=\sum_{i\in\left\{  1,2,\ldots,n\right\}
}=\sum_{i=1}^{n}$ (an equality between summation signs). Now, Proposition
\ref{prop.perm.lehmer.l} yields%
\[
\ell\left(  \sigma\right)  =\ell_{1}\left(  \sigma\right)  +\ell_{2}\left(
\sigma\right)  +\cdots+\ell_{n}\left(  \sigma\right)  =\underbrace{\sum
_{i=1}^{n}}_{=\sum_{i\in\left[  n\right]  }}\ell_{i}\left(  \sigma\right)
=\sum_{i\in\left[  n\right]  }\ell_{i}\left(  \sigma\right)  =\sum
_{i\in\left[  n\right]  }\ell_{\tau\left(  i\right)  }\left(  \sigma\right)
\]
(here, we have substituted $\tau\left(  i\right)  $ for $i$ in the sum, since
the map $\tau:\left[  n\right]  \rightarrow\left[  n\right]  $ is a
bijection). Also, Proposition \ref{prop.perm.lehmer.l} (applied to $\tau$
instead of $\sigma$) yields%
\[
\ell\left(  \tau\right)  =\ell_{1}\left(  \tau\right)  +\ell_{2}\left(
\tau\right)  +\cdots+\ell_{n}\left(  \tau\right)  =\underbrace{\sum_{i=1}^{n}%
}_{=\sum_{i\in\left[  n\right]  }}\ell_{i}\left(  \tau\right)  =\sum
_{i\in\left[  n\right]  }\ell_{i}\left(  \tau\right)  .
\]
Moreover, Proposition \ref{prop.perm.lehmer.l} (applied to $\sigma\circ\tau$
instead of $\sigma$) yields%
\[
\ell\left(  \sigma\circ\tau\right)  =\ell_{1}\left(  \sigma\circ\tau\right)
+\ell_{2}\left(  \sigma\circ\tau\right)  +\cdots+\ell_{n}\left(  \sigma
\circ\tau\right)  =\underbrace{\sum_{i=1}^{n}}_{=\sum_{i\in\left[  n\right]
}}\ell_{i}\left(  \sigma\circ\tau\right)  =\sum_{i\in\left[  n\right]  }%
\ell_{i}\left(  \sigma\circ\tau\right)  .
\]

Now,%
\begin{align}
&  \underbrace{\ell\left(  \sigma\right)  }_{=\sum_{i\in\left[  n\right]
}\ell_{\tau\left(  i\right)  }\left(  \sigma\right)  }+\underbrace{\ell\left(
\tau\right)  }_{=\sum_{i\in\left[  n\right]  }\ell_{i}\left(  \tau\right)
}-\underbrace{\ell\left(  \sigma\circ\tau\right)  }_{=\sum_{i\in\left[
n\right]  }\ell_{i}\left(  \sigma\circ\tau\right)  }\nonumber\\
&  =\sum_{i\in\left[  n\right]  }\ell_{\tau\left(  i\right)  }\left(
\sigma\right)  +\sum_{i\in\left[  n\right]  }\ell_{i}\left(  \tau\right)
-\sum_{i\in\left[  n\right]  }\ell_{i}\left(  \sigma\circ\tau\right)
\nonumber\\
&  =\sum_{i\in\left[  n\right]  }\underbrace{\left(  \ell_{\tau\left(
i\right)  }\left(  \sigma\right)  +\ell_{i}\left(  \tau\right)  -\ell
_{i}\left(  \sigma\circ\tau\right)  \right)  }_{\substack{=\sum_{j\in\left[
n\right]  }\left[  j>i\right]  \left[  \tau\left(  i\right)  >\tau\left(
j\right)  \right]  \left[  \sigma\left(  \tau\left(  j\right)  \right)
>\sigma\left(  \tau\left(  i\right)  \right)  \right]
\\\ \ \ \ \ \ \ \ \ \ +\sum_{j\in\left[  n\right]  }\left[  i>j\right]
\left[  \tau\left(  j\right)  >\tau\left(  i\right)  \right]  \left[
\sigma\left(  \tau\left(  i\right)  \right)  >\sigma\left(  \tau\left(
j\right)  \right)  \right]  \\\text{(by Exercise \ref{exe.perm.lisitau}
\textbf{(a)})}}}\nonumber\\
&  =\sum_{i\in\left[  n\right]  }\left(  \sum_{j\in\left[  n\right]  }\left[
j>i\right]  \left[  \tau\left(  i\right)  >\tau\left(  j\right)  \right]
\left[  \sigma\left(  \tau\left(  j\right)  \right)  >\sigma\left(
\tau\left(  i\right)  \right)  \right]  \right. \nonumber\\
&  \ \ \ \ \ \ \ \ \ \ \ \ \ \ \ \ \ \ \ \ \left.  +\sum_{j\in\left[
n\right]  }\left[  i>j\right]  \left[  \tau\left(  j\right)  >\tau\left(
i\right)  \right]  \left[  \sigma\left(  \tau\left(  i\right)  \right)
>\sigma\left(  \tau\left(  j\right)  \right)  \right]  \right) \nonumber\\
&  =\sum_{i\in\left[  n\right]  }\sum_{j\in\left[  n\right]  }\left[
j>i\right]  \left[  \tau\left(  i\right)  >\tau\left(  j\right)  \right]
\left[  \sigma\left(  \tau\left(  j\right)  \right)  >\sigma\left(
\tau\left(  i\right)  \right)  \right] \nonumber\\
&  \ \ \ \ \ \ \ \ \ \ +\sum_{i\in\left[  n\right]  }\sum_{j\in\left[
n\right]  }\left[  i>j\right]  \left[  \tau\left(  j\right)  >\tau\left(
i\right)  \right]  \left[  \sigma\left(  \tau\left(  i\right)  \right)
>\sigma\left(  \tau\left(  j\right)  \right)  \right]  .
\label{sol.perm.lisitau.b.1}%
\end{align}

\begin{vershort}
But
\begin{align*}
&  \sum_{i\in\left[  n\right]  }\sum_{j\in\left[  n\right]  }\left[
i>j\right]  \left[  \tau\left(  j\right)  >\tau\left(  i\right)  \right]
\left[  \sigma\left(  \tau\left(  i\right)  \right)  >\sigma\left(
\tau\left(  j\right)  \right)  \right] \\
&  =\underbrace{\sum_{j\in\left[  n\right]  }\sum_{i\in\left[  n\right]  }%
}_{=\sum_{i\in\left[  n\right]  }\sum_{j\in\left[  n\right]  }}\left[
j>i\right]  \left[  \tau\left(  i\right)  >\tau\left(  j\right)  \right]
\left[  \sigma\left(  \tau\left(  j\right)  \right)  >\sigma\left(
\tau\left(  i\right)  \right)  \right] \\
&  \ \ \ \ \ \ \ \ \ \ \left(
\begin{array}
[c]{c}%
\text{here, we have renamed the summation indices }i\text{ and }j\\
\text{as }j\text{ and }i\text{, respectively}%
\end{array}
\right) \\
&  =\sum_{i\in\left[  n\right]  }\sum_{j\in\left[  n\right]  }\left[
j>i\right]  \left[  \tau\left(  i\right)  >\tau\left(  j\right)  \right]
\left[  \sigma\left(  \tau\left(  j\right)  \right)  >\sigma\left(
\tau\left(  i\right)  \right)  \right]  .
\end{align*}

\end{vershort}

\begin{verlong}
Next, we observe that
\begin{align*}
&  \sum_{i\in\left[  n\right]  }\sum_{j\in\left[  n\right]  }\left[
i>j\right]  \left[  \tau\left(  j\right)  >\tau\left(  i\right)  \right]
\left[  \sigma\left(  \tau\left(  i\right)  \right)  >\sigma\left(
\tau\left(  j\right)  \right)  \right] \\
&  =\sum_{i\in\left[  n\right]  }\sum_{v\in\left[  n\right]  }\left[
i>v\right]  \left[  \tau\left(  v\right)  >\tau\left(  i\right)  \right]
\left[  \sigma\left(  \tau\left(  i\right)  \right)  >\sigma\left(
\tau\left(  v\right)  \right)  \right] \\
&  \ \ \ \ \ \ \ \ \ \ \left(  \text{here, we have renamed the summation index
}j\text{ as }v\right) \\
&  =\sum_{j\in\left[  n\right]  }\sum_{v\in\left[  n\right]  }\left[
j>v\right]  \left[  \tau\left(  v\right)  >\tau\left(  j\right)  \right]
\left[  \sigma\left(  \tau\left(  j\right)  \right)  >\sigma\left(
\tau\left(  v\right)  \right)  \right] \\
&  \ \ \ \ \ \ \ \ \ \ \left(  \text{here, we have renamed the summation index
}i\text{ as }j\right) \\
&  =\underbrace{\sum_{j\in\left[  n\right]  }\sum_{i\in\left[  n\right]  }%
}_{=\sum_{i\in\left[  n\right]  }\sum_{j\in\left[  n\right]  }}\left[
j>i\right]  \left[  \tau\left(  i\right)  >\tau\left(  j\right)  \right]
\left[  \sigma\left(  \tau\left(  j\right)  \right)  >\sigma\left(
\tau\left(  i\right)  \right)  \right] \\
&  \ \ \ \ \ \ \ \ \ \ \left(  \text{here, we have renamed the summation index
}v\text{ as }i\right) \\
&  =\sum_{i\in\left[  n\right]  }\sum_{j\in\left[  n\right]  }\left[
j>i\right]  \left[  \tau\left(  i\right)  >\tau\left(  j\right)  \right]
\left[  \sigma\left(  \tau\left(  j\right)  \right)  >\sigma\left(
\tau\left(  i\right)  \right)  \right]  .
\end{align*}

\end{verlong}

Hence, (\ref{sol.perm.lisitau.b.1}) becomes%
\begin{align*}
&  \ell\left(  \sigma\right)  +\ell\left(  \tau\right)  -\ell\left(
\sigma\circ\tau\right) \\
&  =\sum_{i\in\left[  n\right]  }\sum_{j\in\left[  n\right]  }\left[
j>i\right]  \left[  \tau\left(  i\right)  >\tau\left(  j\right)  \right]
\left[  \sigma\left(  \tau\left(  j\right)  \right)  >\sigma\left(
\tau\left(  i\right)  \right)  \right] \\
&  \ \ \ \ \ \ \ \ \ \ +\underbrace{\sum_{i\in\left[  n\right]  }\sum
_{j\in\left[  n\right]  }\left[  i>j\right]  \left[  \tau\left(  j\right)
>\tau\left(  i\right)  \right]  \left[  \sigma\left(  \tau\left(  i\right)
\right)  >\sigma\left(  \tau\left(  j\right)  \right)  \right]  }_{=\sum
_{i\in\left[  n\right]  }\sum_{j\in\left[  n\right]  }\left[  j>i\right]
\left[  \tau\left(  i\right)  >\tau\left(  j\right)  \right]  \left[
\sigma\left(  \tau\left(  j\right)  \right)  >\sigma\left(  \tau\left(
i\right)  \right)  \right]  }\\
&  =\sum_{i\in\left[  n\right]  }\sum_{j\in\left[  n\right]  }\left[
j>i\right]  \left[  \tau\left(  i\right)  >\tau\left(  j\right)  \right]
\left[  \sigma\left(  \tau\left(  j\right)  \right)  >\sigma\left(
\tau\left(  i\right)  \right)  \right] \\
&  \ \ \ \ \ \ \ \ \ \ +\sum_{i\in\left[  n\right]  }\sum_{j\in\left[
n\right]  }\left[  j>i\right]  \left[  \tau\left(  i\right)  >\tau\left(
j\right)  \right]  \left[  \sigma\left(  \tau\left(  j\right)  \right)
>\sigma\left(  \tau\left(  i\right)  \right)  \right] \\
&  =2\sum_{i\in\left[  n\right]  }\sum_{j\in\left[  n\right]  }\left[
j>i\right]  \left[  \tau\left(  i\right)  >\tau\left(  j\right)  \right]
\left[  \sigma\left(  \tau\left(  j\right)  \right)  >\sigma\left(
\tau\left(  i\right)  \right)  \right]  .
\end{align*}
This solves Exercise \ref{exe.perm.lisitau} \textbf{(b)}.

\textbf{(c)} \textit{Second solution to Exercise \ref{exe.ps2.2.5}
\textbf{(a)}.} Forget that we fixed $\sigma$ and $\tau$.

Let $k\in\left\{  1,2,\ldots,n-1\right\}  $. Thus, $k+1\in\left\{
2,3,\ldots,n\right\}  \subseteq\left\{  1,2,\ldots,n\right\}  =\left[
n\right]  $. Also, $k\in\left\{  1,2,\ldots,n-1\right\}  \subseteq\left\{
1,2,\ldots,n\right\}  =\left[  n\right]  $.

Recall that $s_{k}$ is the permutation in $S_{n}$ that swaps $k$ with $k+1$
but leaves all other numbers unchanged (by the definition of $s_{k}$). Thus,
we have $s_{k}\left(  k\right)  =k+1$ and $s_{k}\left(  k+1\right)  =k$ and
\[
\left(  s_{k}\left(  j\right)  =j\ \ \ \ \ \ \ \ \ \ \text{for each }%
j\in\left\{  1,2,\ldots,n\right\}  \setminus\left\{  k,k+1\right\}  \right)
.
\]

Let us prove some auxiliary claims:

\begin{statement}
\textit{Claim 1:} We have $\ell\left(  s_{k}\right)  =1$.
\end{statement}

\begin{vershort}
[\textit{Proof of Claim 1:} Let $i\in\left[  n\right]  $. Lemma
\ref{lem.sol.perm.lisitau.1} \textbf{(b)} (applied to $\sigma=s_{k}$) yields%
\[
\sum_{j\in\left[  n\right]  }\left[  i<j\right]  \left[  s_{k}\left(
i\right)  >s_{k}\left(  j\right)  \right]  =\ell_{i}\left(  s_{k}\right)  .
\]
Hence,%
\begin{align}
\ell_{i}\left(  s_{k}\right)   &  =\sum_{j\in\left[  n\right]  }%
\underbrace{\left[  i<j\right]  }_{\substack{=\left[  j>i\right]
\\\text{(since the statement }i<j\\\text{is equivalent to }j>i\text{)}%
}}\left[  s_{k}\left(  i\right)  >s_{k}\left(  j\right)  \right]  =\sum
_{j\in\left[  n\right]  }\underbrace{\left[  j>i\right]  \left[  s_{k}\left(
i\right)  >s_{k}\left(  j\right)  \right]  }_{\substack{=\left[  i=k\right]
\left[  j=k+1\right]  \\\text{(by Lemma \ref{lem.sol.perm.lisitau.sk1}
\textbf{(c)})}}}\nonumber\\
&  =\sum_{j\in\left[  n\right]  }\left[  i=k\right]  \left[  j=k+1\right]
=\left[  i=k\right]  \underbrace{\sum_{j\in\left[  n\right]  }\left[
j=k+1\right]  }_{\substack{=1\\\text{(by Lemma
\ref{lem.sol.perm.lisitau.trivsum},}\\\text{applied to }P=\left[  n\right]
\text{ and }p=k+1\text{)}}}\nonumber\\
&  =\left[  i=k\right]  . \label{sol.perm.lisitau.c.c1.pf.short.4}%
\end{align}

Now, forget that we fixed $i$. We thus have proven
(\ref{sol.perm.lisitau.c.c1.pf.short.4}) for each $i\in\left[  n\right]  $.

Now, Proposition \ref{prop.perm.lehmer.l} (applied to $\sigma=s_{k}$) yields%
\begin{align*}
\ell\left(  s_{k}\right)   &  =\ell_{1}\left(  s_{k}\right)  +\ell_{2}\left(
s_{k}\right)  +\cdots+\ell_{n}\left(  s_{k}\right)  =\sum_{i\in\left[
n\right]  }\underbrace{\ell_{i}\left(  s_{k}\right)  }_{\substack{=\left[
i=k\right]  \\\text{(by (\ref{sol.perm.lisitau.c.c1.pf.short.4}))}}%
}=\sum_{i\in\left[  n\right]  }\left[  i=k\right] \\
&  =\sum_{j\in\left[  n\right]  }\left[  j=k\right]
\ \ \ \ \ \ \ \ \ \ \left(  \text{here, we have renamed the summation index
}i\text{ as }j\right) \\
&  =1\ \ \ \ \ \ \ \ \ \ \left(  \text{by Lemma
\ref{lem.sol.perm.lisitau.trivsum}, applied to }P=\left[  n\right]  \text{ and
}p=k\right)  .
\end{align*}
This proves Claim 1.]
\end{vershort}

\begin{verlong}
[\textit{Proof of Claim 1:} Let $i\in\left[  n\right]  $. Each $j\in\left[
n\right]  $ satisfies%
\begin{equation}
\left[  i<j\right]  =\left[  j>i\right]  \label{sol.perm.lisitau.c.c1.pf.2}%
\end{equation}
\footnote{\textit{Proof of (\ref{sol.perm.lisitau.c.c1.pf.2}):} Let
$j\in\left[  n\right]  $. Clearly, $\left(  j>i\right)  $ and $\left(
i<j\right)  $ are two equivalent logical statements. Hence, Exercise
\ref{exe.iverson-prop} \textbf{(a)} (applied to $\mathcal{A}=\left(
j>i\right)  $ and $\mathcal{B}=\left(  i<j\right)  $) yields $\left[
j>i\right]  =\left[  i<j\right]  $. This proves
(\ref{sol.perm.lisitau.c.c1.pf.2}).}.

Lemma \ref{lem.sol.perm.lisitau.trivsum} (applied to $P=\left[  n\right]  $
and $p=k+1$) yields $\sum_{j\in\left[  n\right]  }\left[  j=k+1\right]  =1$
(since $k+1\in\left[  n\right]  $).

But Lemma \ref{lem.sol.perm.lisitau.1} \textbf{(b)} (applied to $\sigma=s_{k}%
$) yields%
\[
\sum_{j\in\left[  n\right]  }\left[  i<j\right]  \left[  s_{k}\left(
i\right)  >s_{k}\left(  j\right)  \right]  =\ell_{i}\left(  s_{k}\right)  .
\]
Hence,%
\begin{align}
\ell_{i}\left(  s_{k}\right)   &  =\sum_{j\in\left[  n\right]  }%
\underbrace{\left[  i<j\right]  }_{\substack{=\left[  j>i\right]  \\\text{(by
(\ref{sol.perm.lisitau.c.c1.pf.2}))}}}\left[  s_{k}\left(  i\right)
>s_{k}\left(  j\right)  \right]  =\sum_{j\in\left[  n\right]  }%
\underbrace{\left[  j>i\right]  \left[  s_{k}\left(  i\right)  >s_{k}\left(
j\right)  \right]  }_{\substack{=\left[  i=k\right]  \left[  j=k+1\right]
\\\text{(by Lemma \ref{lem.sol.perm.lisitau.sk1} \textbf{(c)})}}}\nonumber\\
&  =\sum_{j\in\left[  n\right]  }\left[  i=k\right]  \left[  j=k+1\right]
=\left[  i=k\right]  \underbrace{\sum_{j\in\left[  n\right]  }\left[
j=k+1\right]  }_{=1}\nonumber\\
&  =\left[  i=k\right]  . \label{sol.perm.lisitau.c.c1.pf.4}%
\end{align}

Now, forget that we fixed $i$. We thus have proven
(\ref{sol.perm.lisitau.c.c1.pf.4}) for each $i\in\left[  n\right]  $.

Lemma \ref{lem.sol.perm.lisitau.trivsum} (applied to $P=\left[  n\right]  $
and $p=k$) yields $\sum_{j\in\left[  n\right]  }\left[  j=k\right]  =1$ (since
$k\in\left[  n\right]  $).

We have $\left[  n\right]  =\left\{  1,2,\ldots,n\right\}  $ and thus
$\sum_{i\in\left[  n\right]  }=\sum_{i\in\left\{  1,2,\ldots,n\right\}  }%
=\sum_{i=1}^{n}$ (an equality between summation signs). Now, Proposition
\ref{prop.perm.lehmer.l} (applied to $\sigma=s_{k}$) yields%
\begin{align*}
\ell\left(  s_{k}\right)   &  =\ell_{1}\left(  s_{k}\right)  +\ell_{2}\left(
s_{k}\right)  +\cdots+\ell_{n}\left(  s_{k}\right)  =\underbrace{\sum
_{i=1}^{n}}_{=\sum_{i\in\left[  n\right]  }}\ell_{i}\left(  s_{k}\right)
=\sum_{i\in\left[  n\right]  }\underbrace{\ell_{i}\left(  s_{k}\right)
}_{\substack{=\left[  i=k\right]  \\\text{(by
(\ref{sol.perm.lisitau.c.c1.pf.4}))}}}=\sum_{i\in\left[  n\right]  }\left[
i=k\right] \\
&  =\sum_{j\in\left[  n\right]  }\left[  j=k\right]
\ \ \ \ \ \ \ \ \ \ \left(  \text{here, we have renamed the summation index
}i\text{ as }j\right) \\
&  =1.
\end{align*}
This proves Claim 1.]
\end{verlong}

\begin{statement}
\textit{Claim 2:} Let $\sigma\in S_{n}$. Then, $\ell\left(  \sigma\circ
s_{k}\right)  =\ell\left(  \sigma\right)  +1-2\left[  \sigma\left(  k\right)
>\sigma\left(  k+1\right)  \right]  $.
\end{statement}

\begin{vershort}
[\textit{Proof of Claim 2:} Every $i\in\left[  n\right]  $ satisfies%
\begin{align*}
&  \sum_{j\in\left[  n\right]  }\underbrace{\left[  j>i\right]  \left[
s_{k}\left(  i\right)  >s_{k}\left(  j\right)  \right]  }_{\substack{=\left[
i=k\right]  \left[  j=k+1\right]  \\\text{(by Lemma
\ref{lem.sol.perm.lisitau.sk1} \textbf{(c)})}}}\left[  \sigma\left(
s_{k}\left(  j\right)  \right)  >\sigma\left(  s_{k}\left(  i\right)  \right)
\right] \\
&  =\sum_{j\in\left[  n\right]  }\left[  i=k\right]  \left[  j=k+1\right]
\left[  \sigma\left(  s_{k}\left(  j\right)  \right)  >\sigma\left(
s_{k}\left(  i\right)  \right)  \right] \\
&  =\left[  i=k\right]  \underbrace{\left[  k+1=k+1\right]  }%
_{\substack{=1\\\text{(since }k+1=k+1\text{)}}}\left[  \sigma\left(
\underbrace{s_{k}\left(  k+1\right)  }_{=k}\right)  >\sigma\left(
s_{k}\left(  i\right)  \right)  \right] \\
&  \ \ \ \ \ \ \ \ \ \ +\sum_{\substack{j\in\left[  n\right]  ;\\j\neq
k+1}}\left[  i=k\right]  \underbrace{\left[  j=k+1\right]  }%
_{\substack{_{\substack{=0}}\\\text{(since }j\neq k+1\text{)}}}\left[
\sigma\left(  s_{k}\left(  j\right)  \right)  >\sigma\left(  s_{k}\left(
i\right)  \right)  \right] \\
&  \ \ \ \ \ \ \ \ \ \ \left(
\begin{array}
[c]{c}%
\text{here, we have split off the addend}\\
\text{for }j=k+1\text{ from the sum (since }k+1\in\left[  n\right]  \text{)}%
\end{array}
\right) \\
&  =\left[  i=k\right]  \left[  \sigma\left(  k\right)  >\sigma\left(
s_{k}\left(  i\right)  \right)  \right]  +\underbrace{\sum_{\substack{j\in
\left[  n\right]  ;\\j\neq k+1}}\left[  i=k\right]  0\left[  \sigma\left(
s_{k}\left(  j\right)  \right)  >\sigma\left(  s_{k}\left(  i\right)  \right)
\right]  }_{=0}\\
&  =\left[  i=k\right]  \left[  \sigma\left(  k\right)  >\sigma\left(
s_{k}\left(  i\right)  \right)  \right]  .
\end{align*}
Summing up these equalities for all $i\in\left[  n\right]  $, we obtain%
\begin{align*}
&  \sum_{i\in\left[  n\right]  }\sum_{j\in\left[  n\right]  }\left[
j>i\right]  \left[  s_{k}\left(  i\right)  >s_{k}\left(  j\right)  \right]
\left[  \sigma\left(  s_{k}\left(  j\right)  \right)  >\sigma\left(
s_{k}\left(  i\right)  \right)  \right] \\
&  =\sum_{i\in\left[  n\right]  }\left[  i=k\right]  \left[  \sigma\left(
k\right)  >\sigma\left(  s_{k}\left(  i\right)  \right)  \right] \\
&  =\underbrace{\left[  k=k\right]  }_{\substack{=1\\\text{(since }%
k=k\text{)}}}\left[  \sigma\left(  k\right)  >\sigma\left(  \underbrace{s_{k}%
\left(  k\right)  }_{=k+1}\right)  \right]  +\sum_{\substack{i\in\left[
n\right]  ;\\i\neq k}}\underbrace{\left[  i=k\right]  }%
_{\substack{=0\\\text{(since }i\neq k\text{)}}}\left[  \sigma\left(  k\right)
>\sigma\left(  s_{k}\left(  i\right)  \right)  \right] \\
&  \ \ \ \ \ \ \ \ \ \ \left(  \text{here, we have split off the addend for
}i=k\text{ from the sum (since }k\in\left[  n\right]  \text{)}\right) \\
&  =\left[  \sigma\left(  k\right)  >\sigma\left(  k+1\right)  \right]
+\underbrace{\sum_{\substack{i\in\left[  n\right]  ;\\i\neq k}}0\left[
\sigma\left(  k\right)  >\sigma\left(  s_{k}\left(  i\right)  \right)
\right]  }_{=0}=\left[  \sigma\left(  k\right)  >\sigma\left(  k+1\right)
\right]  .
\end{align*}
But Exercise \ref{exe.perm.lisitau} \textbf{(b)} (applied to $\tau=s_{k}$)
yields%
\begin{align*}
&  \ell\left(  \sigma\right)  +\ell\left(  s_{k}\right)  -\ell\left(
\sigma\circ s_{k}\right) \\
&  =2\underbrace{\sum_{i\in\left[  n\right]  }\sum_{j\in\left[  n\right]
}\left[  j>i\right]  \left[  s_{k}\left(  i\right)  >s_{k}\left(  j\right)
\right]  \left[  \sigma\left(  s_{k}\left(  j\right)  \right)  >\sigma\left(
s_{k}\left(  i\right)  \right)  \right]  }_{=\left[  \sigma\left(  k\right)
>\sigma\left(  k+1\right)  \right]  }=2\left[  \sigma\left(  k\right)
>\sigma\left(  k+1\right)  \right]  .
\end{align*}
Solving this equation for $\ell\left(  \sigma\circ s_{k}\right)  $, we obtain%
\begin{align*}
\ell\left(  \sigma\circ s_{k}\right)   &  =\ell\left(  \sigma\right)
+\underbrace{\ell\left(  s_{k}\right)  }_{\substack{=1\\\text{(by Claim 1)}%
}}-2\left[  \sigma\left(  k\right)  >\sigma\left(  k+1\right)  \right] \\
&  =\ell\left(  \sigma\right)  +1-2\left[  \sigma\left(  k\right)
>\sigma\left(  k+1\right)  \right]  .
\end{align*}
This proves Claim 2.]
\end{vershort}

\begin{verlong}
[\textit{Proof of Claim 2:} Exercise \ref{exe.perm.lisitau} \textbf{(b)}
(applied to $\tau=s_{k}$) yields%
\begin{align}
&  \ell\left(  \sigma\right)  +\ell\left(  s_{k}\right)  -\ell\left(
\sigma\circ s_{k}\right) \nonumber\\
&  =2\sum_{i\in\left[  n\right]  }\sum_{j\in\left[  n\right]  }\left[
j>i\right]  \left[  s_{k}\left(  i\right)  >s_{k}\left(  j\right)  \right]
\left[  \sigma\left(  s_{k}\left(  j\right)  \right)  >\sigma\left(
s_{k}\left(  i\right)  \right)  \right]  . \label{sol.perm.lisitau.c.c2.pf.1}%
\end{align}

But every $i\in\left[  n\right]  $ satisfies%
\begin{align*}
&  \sum_{j\in\left[  n\right]  }\underbrace{\left[  j>i\right]  \left[
s_{k}\left(  i\right)  >s_{k}\left(  j\right)  \right]  }_{\substack{=\left[
i=k\right]  \left[  j=k+1\right]  \\\text{(by Lemma
\ref{lem.sol.perm.lisitau.sk1} \textbf{(c)})}}}\left[  \sigma\left(
s_{k}\left(  j\right)  \right)  >\sigma\left(  s_{k}\left(  i\right)  \right)
\right] \\
&  =\sum_{j\in\left[  n\right]  }\left[  i=k\right]  \left[  j=k+1\right]
\left[  \sigma\left(  s_{k}\left(  j\right)  \right)  >\sigma\left(
s_{k}\left(  i\right)  \right)  \right] \\
&  =\left[  i=k\right]  \underbrace{\left[  k+1=k+1\right]  }%
_{\substack{=1\\\text{(since }k+1=k+1\text{)}}}\left[  \sigma\left(
\underbrace{s_{k}\left(  k+1\right)  }_{=k}\right)  >\sigma\left(
s_{k}\left(  i\right)  \right)  \right] \\
&  \ \ \ \ \ \ \ \ \ \ +\sum_{\substack{j\in\left[  n\right]  ;\\j\neq
k+1}}\left[  i=k\right]  \underbrace{\left[  j=k+1\right]  }%
_{\substack{=0\\\text{(since we don't have }j=k+1\\\text{(since }j\neq
k+1\text{))}}}\left[  \sigma\left(  s_{k}\left(  j\right)  \right)
>\sigma\left(  s_{k}\left(  i\right)  \right)  \right] \\
&  \ \ \ \ \ \ \ \ \ \ \left(
\begin{array}
[c]{c}%
\text{here, we have split off the addend}\\
\text{for }j=k+1\text{ from the sum (since }k+1\in\left[  n\right]  \text{)}%
\end{array}
\right) \\
&  =\left[  i=k\right]  \left[  \sigma\left(  k\right)  >\sigma\left(
s_{k}\left(  i\right)  \right)  \right]  +\underbrace{\sum_{\substack{j\in
\left[  n\right]  ;\\j\neq k+1}}\left[  i=k\right]  0\left[  \sigma\left(
s_{k}\left(  j\right)  \right)  >\sigma\left(  s_{k}\left(  i\right)  \right)
\right]  }_{=0}\\
&  =\left[  i=k\right]  \left[  \sigma\left(  k\right)  >\sigma\left(
s_{k}\left(  i\right)  \right)  \right]  .
\end{align*}
Summing up these equalities for all $i\in\left[  n\right]  $, we obtain%
\begin{align*}
&  \sum_{i\in\left[  n\right]  }\sum_{j\in\left[  n\right]  }\left[
j>i\right]  \left[  s_{k}\left(  i\right)  >s_{k}\left(  j\right)  \right]
\left[  \sigma\left(  s_{k}\left(  j\right)  \right)  >\sigma\left(
s_{k}\left(  i\right)  \right)  \right] \\
&  =\sum_{i\in\left[  n\right]  }\left[  i=k\right]  \left[  \sigma\left(
k\right)  >\sigma\left(  s_{k}\left(  i\right)  \right)  \right] \\
&  =\underbrace{\left[  k=k\right]  }_{\substack{=1\\\text{(since }%
k=k\text{)}}}\left[  \sigma\left(  k\right)  >\sigma\left(  \underbrace{s_{k}%
\left(  k\right)  }_{=k+1}\right)  \right]  +\sum_{\substack{i\in\left[
n\right]  ;\\i\neq k}}\underbrace{\left[  i=k\right]  }%
_{\substack{=0\\\text{(since we don't have }i=k\\\text{(since }i\neq
k\text{))}}}\left[  \sigma\left(  k\right)  >\sigma\left(  s_{k}\left(
i\right)  \right)  \right] \\
&  \ \ \ \ \ \ \ \ \ \ \left(  \text{here, we have split off the addend for
}i=k\text{ from the sum (since }k\in\left[  n\right]  \text{)}\right) \\
&  =\left[  \sigma\left(  k\right)  >\sigma\left(  k+1\right)  \right]
+\underbrace{\sum_{\substack{i\in\left[  n\right]  ;\\i\neq k}}0\left[
\sigma\left(  k\right)  >\sigma\left(  s_{k}\left(  i\right)  \right)
\right]  }_{=0}=\left[  \sigma\left(  k\right)  >\sigma\left(  k+1\right)
\right]  .
\end{align*}
Hence, (\ref{sol.perm.lisitau.c.c2.pf.1}) becomes%
\begin{align*}
&  \ell\left(  \sigma\right)  +\ell\left(  s_{k}\right)  -\ell\left(
\sigma\circ s_{k}\right) \\
&  =2\underbrace{\sum_{i\in\left[  n\right]  }\sum_{j\in\left[  n\right]
}\left[  j>i\right]  \left[  s_{k}\left(  i\right)  >s_{k}\left(  j\right)
\right]  \left[  \sigma\left(  s_{k}\left(  j\right)  \right)  >\sigma\left(
s_{k}\left(  i\right)  \right)  \right]  }_{=\left[  \sigma\left(  k\right)
>\sigma\left(  k+1\right)  \right]  }=2\left[  \sigma\left(  k\right)
>\sigma\left(  k+1\right)  \right]  .
\end{align*}
Solving this equation for $\ell\left(  \sigma\circ s_{k}\right)  $, we obtain%
\begin{align*}
\ell\left(  \sigma\circ s_{k}\right)   &  =\ell\left(  \sigma\right)
+\underbrace{\ell\left(  s_{k}\right)  }_{\substack{=1\\\text{(by Claim 1)}%
}}-2\left[  \sigma\left(  k\right)  >\sigma\left(  k+1\right)  \right] \\
&  =\ell\left(  \sigma\right)  +1-2\left[  \sigma\left(  k\right)
>\sigma\left(  k+1\right)  \right]  .
\end{align*}
This proves Claim 2.]
\end{verlong}

\begin{statement}
\textit{Claim 3:} Let $\sigma\in S_{n}$. Then,%
\[
\ell\left(  \sigma\circ s_{k}\right)  =%
\begin{cases}
\ell\left(  \sigma\right)  +1, & \text{if }\sigma\left(  k\right)
<\sigma\left(  k+1\right)  ;\\
\ell\left(  \sigma\right)  -1, & \text{if }\sigma\left(  k\right)
>\sigma\left(  k+1\right)
\end{cases}
.
\]

\end{statement}

\begin{vershort}
[\textit{Proof of Claim 3:} We have $\sigma\in S_{n}$. In other words,
$\sigma$ is a permutation of $\left\{  1,2,\ldots,n\right\}  $. Thus, $\sigma$
is a bijective map, and therefore an injective map. Hence, from $k\neq k+1$,
we obtain $\sigma\left(  k\right)  \neq\sigma\left(  k+1\right)  $. Thus, the
statement $\left(  \sigma\left(  k\right)  \leq\sigma\left(  k+1\right)
\right)  $ is equivalent to the statement $\left(  \sigma\left(  k\right)
<\sigma\left(  k+1\right)  \right)  $. Moreover, Claim 2 yields%
\begin{align*}
\ell\left(  \sigma\circ s_{k}\right)   &  =\ell\left(  \sigma\right)
+1-2\underbrace{\left[  \sigma\left(  k\right)  >\sigma\left(  k+1\right)
\right]  }_{\substack{=%
\begin{cases}
0, & \text{if }\sigma\left(  k\right)  \leq\sigma\left(  k+1\right)  ;\\
1, & \text{if }\sigma\left(  k\right)  >\sigma\left(  k+1\right)
\end{cases}
\\=%
\begin{cases}
0, & \text{if }\sigma\left(  k\right)  <\sigma\left(  k+1\right)  ;\\
1, & \text{if }\sigma\left(  k\right)  >\sigma\left(  k+1\right)
\end{cases}
\\\text{(since the statement }\left(  \sigma\left(  k\right)  \leq
\sigma\left(  k+1\right)  \right)  \text{ is equivalent}\\\text{to the
statement }\left(  \sigma\left(  k\right)  <\sigma\left(  k+1\right)  \right)
\text{)}}}\\
&  =\ell\left(  \sigma\right)  +1-2%
\begin{cases}
0, & \text{if }\sigma\left(  k\right)  <\sigma\left(  k+1\right)  ;\\
1, & \text{if }\sigma\left(  k\right)  >\sigma\left(  k+1\right)
\end{cases}
\\
&  =%
\begin{cases}
\ell\left(  \sigma\right)  +1-2\cdot0, & \text{if }\sigma\left(  k\right)
<\sigma\left(  k+1\right)  ;\\
\ell\left(  \sigma\right)  +1-2\cdot1, & \text{if }\sigma\left(  k\right)
>\sigma\left(  k+1\right)
\end{cases}
\\
&  =%
\begin{cases}
\ell\left(  \sigma\right)  +1, & \text{if }\sigma\left(  k\right)
<\sigma\left(  k+1\right)  ;\\
\ell\left(  \sigma\right)  -1, & \text{if }\sigma\left(  k\right)
>\sigma\left(  k+1\right)
\end{cases}
.
\end{align*}
This proves Claim 3.]
\end{vershort}

\begin{verlong}
[\textit{Proof of Claim 3:} We have $\sigma\in S_{n}$. In other words,
$\sigma$ is a permutation of $\left\{  1,2,\ldots,n\right\}  $ (since $S_{n}$
is the set of all permutations of $\left\{  1,2,\ldots,n\right\}  $). In other
words, $\sigma$ is a bijection $\left\{  1,2,\ldots,n\right\}  \rightarrow
\left\{  1,2,\ldots,n\right\}  $. Hence, the map $\sigma$ is bijective,
therefore injective. Hence, if we had $\sigma\left(  k\right)  =\sigma\left(
k+1\right)  $, then we would have $k=k+1$, which would contradict $k\neq k+1$.
Thus, we don't have $\sigma\left(  k\right)  =\sigma\left(  k+1\right)  $.
Hence, we have $\sigma\left(  k\right)  \neq\sigma\left(  k+1\right)  $. Thus,
we have either $\sigma\left(  k\right)  <\sigma\left(  k+1\right)  $ or
$\sigma\left(  k\right)  >\sigma\left(  k+1\right)  $. In other words, we are
in one of the following two cases:

\textit{Case 1:} We have $\sigma\left(  k\right)  <\sigma\left(  k+1\right)  $.

\textit{Case 2:} We have $\sigma\left(  k\right)  >\sigma\left(  k+1\right)  $.

Let us first consider Case 1. In this case, we have $\sigma\left(  k\right)
<\sigma\left(  k+1\right)  $. Hence, we don't have $\sigma\left(  k\right)
>\sigma\left(  k+1\right)  $. Therefore, $\left[  \sigma\left(  k\right)
>\sigma\left(  k+1\right)  \right]  =0$. Now, Claim 2 yields
\[
\ell\left(  \sigma\circ s_{k}\right)  =\ell\left(  \sigma\right)
+1-2\underbrace{\left[  \sigma\left(  k\right)  >\sigma\left(  k+1\right)
\right]  }_{=0}=\ell\left(  \sigma\right)  +1-2\cdot0=\ell\left(
\sigma\right)  +1.
\]
Comparing this with%
\[%
\begin{cases}
\ell\left(  \sigma\right)  +1, & \text{if }\sigma\left(  k\right)
<\sigma\left(  k+1\right)  ;\\
\ell\left(  \sigma\right)  -1, & \text{if }\sigma\left(  k\right)
>\sigma\left(  k+1\right)
\end{cases}
=\ell\left(  \sigma\right)  +1\ \ \ \ \ \ \ \ \ \ \left(  \text{since }%
\sigma\left(  k\right)  <\sigma\left(  k+1\right)  \right)  ,
\]
we obtain $\ell\left(  \sigma\circ s_{k}\right)  =%
\begin{cases}
\ell\left(  \sigma\right)  +1, & \text{if }\sigma\left(  k\right)
<\sigma\left(  k+1\right)  ;\\
\ell\left(  \sigma\right)  -1, & \text{if }\sigma\left(  k\right)
>\sigma\left(  k+1\right)
\end{cases}
$. Thus, Claim 3 is proven in Case 1.

Let us next consider Case 2. In this case, we have $\sigma\left(  k\right)
>\sigma\left(  k+1\right)  $. Hence, $\left[  \sigma\left(  k\right)
>\sigma\left(  k+1\right)  \right]  =1$. Now, Claim 2 yields
\[
\ell\left(  \sigma\circ s_{k}\right)  =\ell\left(  \sigma\right)
+1-2\underbrace{\left[  \sigma\left(  k\right)  >\sigma\left(  k+1\right)
\right]  }_{=1}=\ell\left(  \sigma\right)  +1-2\cdot1=\ell\left(
\sigma\right)  -1.
\]
Comparing this with%
\[%
\begin{cases}
\ell\left(  \sigma\right)  +1, & \text{if }\sigma\left(  k\right)
<\sigma\left(  k+1\right)  ;\\
\ell\left(  \sigma\right)  -1, & \text{if }\sigma\left(  k\right)
>\sigma\left(  k+1\right)
\end{cases}
=\ell\left(  \sigma\right)  -1\ \ \ \ \ \ \ \ \ \ \left(  \text{since }%
\sigma\left(  k\right)  >\sigma\left(  k+1\right)  \right)  ,
\]
we obtain $\ell\left(  \sigma\circ s_{k}\right)  =%
\begin{cases}
\ell\left(  \sigma\right)  +1, & \text{if }\sigma\left(  k\right)
<\sigma\left(  k+1\right)  ;\\
\ell\left(  \sigma\right)  -1, & \text{if }\sigma\left(  k\right)
>\sigma\left(  k+1\right)
\end{cases}
$. Thus, Claim 3 is proven in Case 2.

We have now proven Claim 3 in both Cases 1 and 2. Since these two Cases cover
all possibilities, we thus conclude that Claim 3 always holds.]
\end{verlong}

\begin{statement}
\textit{Claim 4:} Let $\sigma\in S_{n}$. Then,%
\[
\ell\left(  s_{k}\circ\sigma\right)  =%
\begin{cases}
\ell\left(  \sigma\right)  +1, & \text{if }\sigma^{-1}\left(  k\right)
<\sigma^{-1}\left(  k+1\right)  ;\\
\ell\left(  \sigma\right)  -1, & \text{if }\sigma^{-1}\left(  k\right)
>\sigma^{-1}\left(  k+1\right)
\end{cases}
.
\]

\end{statement}

[\textit{Proof of Claim 4:} Recall that $s_{i}^{2}=\operatorname*{id}$ for
each $i\in\left\{  1,2,\ldots,n-1\right\}  $. Applying this to $i=k$, we
obtain $s_{k}^{2}=\operatorname*{id}$; thus, $s_{k}\circ s_{k}=s_{k}%
^{2}=\operatorname*{id}$ and therefore $s_{k}^{-1}=s_{k}$.

Let us recall that $\left(  \alpha\circ\beta\right)  ^{-1}=\beta^{-1}%
\circ\alpha^{-1}$ for any two permutations $\alpha$ and $\beta$ in $S_{n}$.
Applying this to $\alpha=\sigma^{-1}$ and $\beta=s_{k}$, we obtain $\left(
\sigma^{-1}\circ s_{k}\right)  ^{-1}=\underbrace{s_{k}^{-1}}_{=s_{k}}%
\circ\underbrace{\left(  \sigma^{-1}\right)  ^{-1}}_{=\sigma}=s_{k}\circ
\sigma$. But Exercise \ref{exe.ps2.2.5} \textbf{(f)} yields $\ell\left(
\sigma\right)  =\ell\left(  \sigma^{-1}\right)  $. Also, Exercise
\ref{exe.ps2.2.5} \textbf{(f)} (applied to $\sigma^{-1}\circ s_{k}$ instead of
$\sigma$) yields $\ell\left(  \sigma^{-1}\circ s_{k}\right)  =\ell\left(
\underbrace{\left(  \sigma^{-1}\circ s_{k}\right)  ^{-1}}_{=s_{k}\circ\sigma
}\right)  =\ell\left(  s_{k}\circ\sigma\right)  $. But Claim 3 (applied to
$\sigma^{-1}$ instead of $\sigma$) yields%
\[
\ell\left(  \sigma^{-1}\circ s_{k}\right)  =%
\begin{cases}
\ell\left(  \sigma^{-1}\right)  +1, & \text{if }\sigma^{-1}\left(  k\right)
<\sigma^{-1}\left(  k+1\right)  ;\\
\ell\left(  \sigma^{-1}\right)  -1, & \text{if }\sigma^{-1}\left(  k\right)
>\sigma^{-1}\left(  k+1\right)
\end{cases}
.
\]
Since $\ell\left(  \sigma^{-1}\circ s_{k}\right)  =\ell\left(  s_{k}%
\circ\sigma\right)  $ and $\ell\left(  \sigma^{-1}\right)  =\ell\left(
\sigma\right)  $, this equality rewrites as follows:%
\[
\ell\left(  s_{k}\circ\sigma\right)  =%
\begin{cases}
\ell\left(  \sigma\right)  +1, & \text{if }\sigma^{-1}\left(  k\right)
<\sigma^{-1}\left(  k+1\right)  ;\\
\ell\left(  \sigma\right)  -1, & \text{if }\sigma^{-1}\left(  k\right)
>\sigma^{-1}\left(  k+1\right)
\end{cases}
.
\]
This proves Claim 4.]

Combining Claim 3 with Claim 4, we obtain the exact statement of Exercise
\ref{exe.ps2.2.5} \textbf{(a)}. Thus, Exercise \ref{exe.ps2.2.5} \textbf{(a)}
is solved again. This solves Exercise \ref{exe.perm.lisitau} \textbf{(c)}.

\begin{vershort}
\textbf{(d)} \textit{Second solution to Exercise \ref{exe.ps2.2.5}
\textbf{(b)}.} Let $\sigma\in S_{n}$ and $\tau\in S_{n}$. Exercise
\ref{exe.perm.lisitau} \textbf{(b)} yields%
\begin{align*}
\ell\left(  \sigma\right)  +\ell\left(  \tau\right)  -\ell\left(  \sigma
\circ\tau\right)   &  =2\underbrace{\sum_{i\in\left[  n\right]  }\sum
_{j\in\left[  n\right]  }\left[  j>i\right]  \left[  \tau\left(  i\right)
>\tau\left(  j\right)  \right]  \left[  \sigma\left(  \tau\left(  j\right)
\right)  >\sigma\left(  \tau\left(  i\right)  \right)  \right]  }%
_{\substack{\text{This is an integer}\\\text{(since the truth values }\left[
j>i\right]  \text{, }\left[  \tau\left(  i\right)  >\tau\left(  j\right)
\right]  \text{, }\left[  \sigma\left(  \tau\left(  j\right)  \right)
>\sigma\left(  \tau\left(  i\right)  \right)  \right]  \text{ are integers)}%
}}\\
&  \equiv0\operatorname{mod}2.
\end{align*}
In other words, $\ell\left(  \sigma\circ\tau\right)  \equiv\ell\left(
\sigma\right)  +\ell\left(  \tau\right)  \operatorname{mod}2$. Thus, Exercise
\ref{exe.ps2.2.5} \textbf{(b)} is solved again. This solves Exercise
\ref{exe.perm.lisitau} \textbf{(d)}.
\end{vershort}

\begin{verlong}
\textbf{(d)} \textit{Second solution to Exercise \ref{exe.ps2.2.5}
\textbf{(b)}.} Forget that we fixed $\sigma$ and $\tau$.

Let $\sigma$ and $\tau$ be two permutations in $S_{n}$. Let%
\begin{equation}
z=\sum_{i\in\left[  n\right]  }\sum_{j\in\left[  n\right]  }\left[
j>i\right]  \left[  \tau\left(  i\right)  >\tau\left(  j\right)  \right]
\left[  \sigma\left(  \tau\left(  j\right)  \right)  >\sigma\left(
\tau\left(  i\right)  \right)  \right]  . \label{sol.perm.lisitau.d.1}%
\end{equation}
Then, $z$ is an integer (because the truth values $\left[  j>i\right]  $,
$\left[  \tau\left(  i\right)  >\tau\left(  j\right)  \right]  $ and $\left[
\sigma\left(  \tau\left(  j\right)  \right)  >\sigma\left(  \tau\left(
i\right)  \right)  \right]  $ are integers for all $i\in\left[  n\right]  $
and $j\in\left[  n\right]  $).

Exercise \ref{exe.perm.lisitau} \textbf{(b)} yields%
\begin{align*}
\ell\left(  \sigma\right)  +\ell\left(  \tau\right)  -\ell\left(  \sigma
\circ\tau\right)   &  =2\underbrace{\sum_{i\in\left[  n\right]  }\sum
_{j\in\left[  n\right]  }\left[  j>i\right]  \left[  \tau\left(  i\right)
>\tau\left(  j\right)  \right]  \left[  \sigma\left(  \tau\left(  j\right)
\right)  >\sigma\left(  \tau\left(  i\right)  \right)  \right]  }%
_{\substack{=z\\\text{(by (\ref{sol.perm.lisitau.d.1}))}}}=2z\\
&  \equiv0\operatorname{mod}2
\end{align*}
(since $z$ is an integer). In other words, $\ell\left(  \sigma\circ
\tau\right)  \equiv\ell\left(  \sigma\right)  +\ell\left(  \tau\right)
\operatorname{mod}2$. Thus, Exercise \ref{exe.ps2.2.5} \textbf{(b)} is solved
again. This solves Exercise \ref{exe.perm.lisitau} \textbf{(d)}.
\end{verlong}

\begin{vershort}
\textbf{(e)} \textit{Second solution to Exercise \ref{exe.ps2.2.5}
\textbf{(c)}.} Let $\sigma\in S_{n}$ and $\tau\in S_{n}$. Exercise
\ref{exe.perm.lisitau} \textbf{(b)} yields%
\[
\ell\left(  \sigma\right)  +\ell\left(  \tau\right)  -\ell\left(  \sigma
\circ\tau\right)  =2\sum_{i\in\left[  n\right]  }\sum_{j\in\left[  n\right]
}\underbrace{\left[  j>i\right]  }_{\geq0}\underbrace{\left[  \tau\left(
i\right)  >\tau\left(  j\right)  \right]  }_{\geq0}\underbrace{\left[
\sigma\left(  \tau\left(  j\right)  \right)  >\sigma\left(  \tau\left(
i\right)  \right)  \right]  }_{\geq0}\geq0.
\]
In other words, $\ell\left(  \sigma\circ\tau\right)  \leq\ell\left(
\sigma\right)  +\ell\left(  \tau\right)  $. Thus, Exercise \ref{exe.ps2.2.5}
\textbf{(c)} is solved again. This solves Exercise \ref{exe.perm.lisitau}
\textbf{(e)}. \qedhere

\end{vershort}

\begin{verlong}
\textbf{(e)} \textit{Second solution to Exercise \ref{exe.ps2.2.5}
\textbf{(c)}.} Forget that we fixed $\sigma$ and $\tau$.

If $\mathcal{A}$ is any logical statement, then%
\begin{equation}
\text{the truth value }\left[  \mathcal{A}\right]  \text{ is nonnegative}
\label{sol.perm.lisitau.e.nn}%
\end{equation}
(since $\left[  \mathcal{A}\right]  \in\left\{  0,1\right\}  $).

Let $\sigma$ and $\tau$ be two permutations in $S_{n}$. Let%
\begin{equation}
z=\sum_{i\in\left[  n\right]  }\sum_{j\in\left[  n\right]  }\left[
j>i\right]  \left[  \tau\left(  i\right)  >\tau\left(  j\right)  \right]
\left[  \sigma\left(  \tau\left(  j\right)  \right)  >\sigma\left(
\tau\left(  i\right)  \right)  \right]  . \label{sol.perm.lisitau.e.1}%
\end{equation}
Then, $z$ is nonnegative (because the truth values $\left[  j>i\right]  $,
$\left[  \tau\left(  i\right)  >\tau\left(  j\right)  \right]  $ and $\left[
\sigma\left(  \tau\left(  j\right)  \right)  >\sigma\left(  \tau\left(
i\right)  \right)  \right]  $ are nonnegative for all $i\in\left[  n\right]  $
and $j\in\left[  n\right]  $ (by (\ref{sol.perm.lisitau.e.nn}))). In other
words, $z\geq0$.

Exercise \ref{exe.perm.lisitau} \textbf{(b)} yields%
\begin{align*}
\ell\left(  \sigma\right)  +\ell\left(  \tau\right)  -\ell\left(  \sigma
\circ\tau\right)   &  =2\underbrace{\sum_{i\in\left[  n\right]  }\sum
_{j\in\left[  n\right]  }\left[  j>i\right]  \left[  \tau\left(  i\right)
>\tau\left(  j\right)  \right]  \left[  \sigma\left(  \tau\left(  j\right)
\right)  >\sigma\left(  \tau\left(  i\right)  \right)  \right]  }%
_{\substack{=z\\\text{(by (\ref{sol.perm.lisitau.e.1}))}}}\\
&  =2z\geq0
\end{align*}
(since $z\geq0$). In other words, $\ell\left(  \sigma\circ\tau\right)
\leq\ell\left(  \sigma\right)  +\ell\left(  \tau\right)  $. Thus, Exercise
\ref{exe.ps2.2.5} \textbf{(c)} is solved again. This solves Exercise
\ref{exe.perm.lisitau} \textbf{(e)}.
\end{verlong}
\end{proof}

\subsection{\label{sect.sol.perm.lisitij}Solution to Exercise
\ref{exe.perm.lisitij}}

Throughout Section \ref{sect.sol.perm.lisitij}, we shall use the same
notations that were in use throughout Section \ref{sect.perm.lehmer}. We shall
furthermore use the notation from Definition \ref{def.iverson}. We need
several lemmas to prepare for the solution of Exercise \ref{exe.perm.lisitij}.
We start with a lemma that restates the definition of $\ell_{i}\left(
\sigma\right)  $:

\begin{lemma}
\label{lem.sol.perm.lisitij.li=}Let $n\in\mathbb{N}$. Let $\sigma\in S_{n}$.
Let $i\in\left[  n\right]  $. Then,%
\[
\ell_{i}\left(  \sigma\right)  =\sum_{j=i+1}^{n}\left[  \sigma\left(
i\right)  >\sigma\left(  j\right)  \right]  .
\]

\end{lemma}

\begin{proof}
[Proof of Lemma \ref{lem.sol.perm.lisitij.li=}.]Recall that $\ell_{i}\left(
\sigma\right)  $ is the number of all $j\in\left\{  i+1,i+2,\ldots,n\right\}
$ such that $\sigma\left(  i\right)  >\sigma\left(  j\right)  $ (by the
definition of $\ell_{i}\left(  \sigma\right)  $). In other words,%
\[
\ell_{i}\left(  \sigma\right)  =\left\vert \left\{  j\in\left\{
i+1,i+2,\ldots,n\right\}  \ \mid\ \sigma\left(  i\right)  >\sigma\left(
j\right)  \right\}  \right\vert .
\]

Comparing this with%
\begin{align*}
&  \underbrace{\sum_{j=i+1}^{n}}_{=\sum_{j\in\left\{  i+1,i+2,\ldots
,n\right\}  }}\left[  \sigma\left(  i\right)  >\sigma\left(  j\right)  \right]
\\
&  =\sum_{j\in\left\{  i+1,i+2,\ldots,n\right\}  }\left[  \sigma\left(
i\right)  >\sigma\left(  j\right)  \right] \\
&  =\sum_{\substack{j\in\left\{  i+1,i+2,\ldots,n\right\}  ;\\\sigma\left(
i\right)  >\sigma\left(  j\right)  }}\underbrace{\left[  \sigma\left(
i\right)  >\sigma\left(  j\right)  \right]  }_{\substack{=1\\\text{(since
}\sigma\left(  i\right)  >\sigma\left(  j\right)  \text{)}}}+\sum
_{\substack{j\in\left\{  i+1,i+2,\ldots,n\right\}  ;\\\text{not }\sigma\left(
i\right)  >\sigma\left(  j\right)  }}\underbrace{\left[  \sigma\left(
i\right)  >\sigma\left(  j\right)  \right]  }_{\substack{=0\\\text{(since we
don't have }\sigma\left(  i\right)  >\sigma\left(  j\right)  \text{)}}}\\
&  \ \ \ \ \ \ \ \ \ \ \left(
\begin{array}
[c]{c}%
\text{since each }j\in\left\{  i+1,i+2,\ldots,n\right\}  \text{ satisfies
either }\left(  \sigma\left(  i\right)  >\sigma\left(  j\right)  \right) \\
\text{or }\left(  \text{not }\sigma\left(  i\right)  >\sigma\left(  j\right)
\right)  \text{ (but not both)}%
\end{array}
\right) \\
&  =\sum_{\substack{j\in\left\{  i+1,i+2,\ldots,n\right\}  ;\\\sigma\left(
i\right)  >\sigma\left(  j\right)  }}1+\underbrace{\sum_{\substack{j\in
\left\{  i+1,i+2,\ldots,n\right\}  ;\\\text{not }\sigma\left(  i\right)
>\sigma\left(  j\right)  }}0}_{=0}=\sum_{\substack{j\in\left\{  i+1,i+2,\ldots
,n\right\}  ;\\\sigma\left(  i\right)  >\sigma\left(  j\right)  }}1\\
&  =\left\vert \left\{  j\in\left\{  i+1,i+2,\ldots,n\right\}  \ \mid
\ \sigma\left(  i\right)  >\sigma\left(  j\right)  \right\}  \right\vert
\cdot1\\
&  =\left\vert \left\{  j\in\left\{  i+1,i+2,\ldots,n\right\}  \ \mid
\ \sigma\left(  i\right)  >\sigma\left(  j\right)  \right\}  \right\vert ,
\end{align*}
we obtain $\ell_{i}\left(  \sigma\right)  =\sum_{j=i+1}^{n}\left[
\sigma\left(  i\right)  >\sigma\left(  j\right)  \right]  $. This proves Lemma
\ref{lem.sol.perm.lisitij.li=}.
\end{proof}

\begin{lemma}
\label{lem.sol.perm.lisitij.uiv}Let $n\in\mathbb{N}$. Let $\sigma\in S_{n}$.
Let $u\in\left[  n\right]  $ and $v\in\left[  n\right]  $ be such that $u<v$.
Let $i\in\left\{  u+1,u+2,\ldots,v-1\right\}  $. Then,%
\begin{align*}
&  \left[  \sigma\left(  u\right)  >\sigma\left(  i\right)  \right]  -\left[
\sigma\left(  v\right)  >\sigma\left(  i\right)  \right]  +\left[
\sigma\left(  i\right)  >\sigma\left(  v\right)  \right]  -\left[
\sigma\left(  i\right)  >\sigma\left(  u\right)  \right] \\
&  =2\left[  \sigma\left(  u\right)  >\sigma\left(  i\right)  >\sigma\left(
v\right)  \right]  -2\left[  \sigma\left(  v\right)  >\sigma\left(  i\right)
>\sigma\left(  u\right)  \right]  .
\end{align*}

\end{lemma}

\begin{vershort}
\begin{proof}
[Proof of Lemma \ref{lem.sol.perm.lisitij.uiv}.]We know that $\sigma$ is a
permutation (since $\sigma\in S_{n}$), thus a bijective map, thus an injective
map. From $i\in\left\{  u+1,u+2,\ldots,v-1\right\}  $, we obtain $u+1\leq
i\leq v-1$. From $u<u+1\leq i$, we obtain $u\neq i$ and thus $\sigma\left(
u\right)  \neq\sigma\left(  i\right)  $ (since $\sigma$ is injective). From
$i\leq v-1<v$, we obtain $i\neq v$ and thus $\sigma\left(  i\right)
\neq\sigma\left(  v\right)  $ (since $\sigma$ is injective).

From $\sigma\left(  u\right)  \neq\sigma\left(  i\right)  $, we conclude that
either $\sigma\left(  u\right)  <\sigma\left(  i\right)  $ or $\sigma\left(
u\right)  >\sigma\left(  i\right)  $. From $\sigma\left(  i\right)  \neq
\sigma\left(  v\right)  $, we conclude that either $\sigma\left(  i\right)
<\sigma\left(  v\right)  $ or $\sigma\left(  i\right)  >\sigma\left(
v\right)  $. Combining the previous two sentences, we conclude that we must be
in one of the following four cases:

\textit{Case 1:} We have $\sigma\left(  u\right)  <\sigma\left(  i\right)  $
and $\sigma\left(  i\right)  <\sigma\left(  v\right)  $.

\textit{Case 2:} We have $\sigma\left(  u\right)  <\sigma\left(  i\right)  $
and $\sigma\left(  i\right)  >\sigma\left(  v\right)  $.

\textit{Case 3:} We have $\sigma\left(  u\right)  >\sigma\left(  i\right)  $
and $\sigma\left(  i\right)  <\sigma\left(  v\right)  $.

\textit{Case 4:} We have $\sigma\left(  u\right)  >\sigma\left(  i\right)  $
and $\sigma\left(  i\right)  >\sigma\left(  v\right)  $.

We can now verify the claim of Lemma \ref{lem.sol.perm.lisitij.uiv} in each of
these four Cases by inspection:

\begin{itemize}
\item In Case 1, the equality claimed in Lemma \ref{lem.sol.perm.lisitij.uiv}
boils down to $0-1+0-1=2\cdot0-2\cdot1$, which is correct.

\item In Case 2, the equality claimed in Lemma \ref{lem.sol.perm.lisitij.uiv}
boils down to $0-0+1-1=2\cdot0-2\cdot0$, which is correct.

\item In Case 3, the equality claimed in Lemma \ref{lem.sol.perm.lisitij.uiv}
boils down to $1-1+0-0=2\cdot0-2\cdot0$, which is correct.

\item In Case 4, the equality claimed in Lemma \ref{lem.sol.perm.lisitij.uiv}
boils down to $1-0+1-0=2\cdot1-2\cdot0$, which is correct.
\end{itemize}

Thus, the equality claimed in Lemma \ref{lem.sol.perm.lisitij.uiv} holds in
all four Cases 1, 2, 3 and 4. Hence, Lemma \ref{lem.sol.perm.lisitij.uiv} is proven.
\end{proof}
\end{vershort}

\begin{verlong}
\begin{proof}
[Proof of Lemma \ref{lem.sol.perm.lisitij.uiv}.]We have $\sigma\in S_{n}$. In
other words, $\sigma$ is a permutation of $\left\{  1,2,\ldots,n\right\}  $
(since $S_{n}$ is the set of all permutations of $\left\{  1,2,\ldots
,n\right\}  $). In other words, $\sigma$ is a bijection $\left\{
1,2,\ldots,n\right\}  \rightarrow\left\{  1,2,\ldots,n\right\}  $. Hence, the
map $\sigma$ is bijective, therefore injective.

We have $u\in\left[  n\right]  =\left\{  1,2,\ldots,n\right\}  $ (by the
definition of $\left[  n\right]  $), thus $1\leq u\leq n$. The same argument
(applied to $v$ instead of $u$) yields $1\leq v\leq n$. But $i\in\left\{
u+1,u+2,\ldots,v-1\right\}  $, so that $u+1\leq i\leq v-1$. Hence, $1\leq
u\leq u+1\leq i$. Combining this with $i\leq v-1\leq v\leq n$, we obtain
$1\leq i\leq n$, so that $i\in\left\{  1,2,\ldots,n\right\}  $. Hence,
$\sigma\left(  i\right)  $ is well-defined.

Also, $u<u+1\leq i$, so that $u\neq i$. If we had $\sigma\left(  u\right)
=\sigma\left(  i\right)  $, then we would have $u=i$ (since $\sigma$ is
injective), which would contradict $u\neq i$. Hence, we cannot have
$\sigma\left(  u\right)  =\sigma\left(  i\right)  $. Thus, we have
$\sigma\left(  u\right)  \neq\sigma\left(  i\right)  $.

Also, $i\leq v-1<v$, so that $i\neq v$. If we had $\sigma\left(  i\right)
=\sigma\left(  v\right)  $, then we would have $i=v$ (since $\sigma$ is
injective), which would contradict $i\neq v$. Hence, we cannot have
$\sigma\left(  i\right)  =\sigma\left(  v\right)  $. Thus, we have
$\sigma\left(  i\right)  \neq\sigma\left(  v\right)  $.

We have $\sigma\left(  u\right)  \neq\sigma\left(  i\right)  $. In other
words, either $\sigma\left(  u\right)  <\sigma\left(  i\right)  $ or
$\sigma\left(  u\right)  >\sigma\left(  i\right)  $. Thus, we are in one of
the following two cases:

\textit{Case 1:} We have $\sigma\left(  u\right)  <\sigma\left(  i\right)  $.

\textit{Case 2:} We have $\sigma\left(  u\right)  >\sigma\left(  i\right)  $.

Let us first consider Case 1. In this case, we have $\sigma\left(  u\right)
<\sigma\left(  i\right)  $. In other words, $\sigma\left(  i\right)
>\sigma\left(  u\right)  $. Hence, $\left[  \sigma\left(  i\right)
>\sigma\left(  u\right)  \right]  =1$. Also, we don't have $\sigma\left(
u\right)  >\sigma\left(  i\right)  $ (since we have $\sigma\left(  u\right)
<\sigma\left(  i\right)  $). Hence, we have $\left[  \sigma\left(  u\right)
>\sigma\left(  i\right)  \right]  =0$. Also, we don't have $\sigma\left(
u\right)  >\sigma\left(  i\right)  >\sigma\left(  v\right)  $ (because we
don't have $\sigma\left(  u\right)  >\sigma\left(  i\right)  $). Thus,
$\left[  \sigma\left(  u\right)  >\sigma\left(  i\right)  >\sigma\left(
v\right)  \right]  =0$.

We have $\sigma\left(  i\right)  \neq\sigma\left(  v\right)  $. In other
words, either $\sigma\left(  i\right)  <\sigma\left(  v\right)  $ or
$\sigma\left(  i\right)  >\sigma\left(  v\right)  $. Thus, we are in one of
the following two subcases:

\textit{Subcase 1.1:} We have $\sigma\left(  i\right)  <\sigma\left(
v\right)  $.

\textit{Subcase 1.2:} We have $\sigma\left(  i\right)  >\sigma\left(
v\right)  $.

Let us first consider Subcase 1.1. In this subcase, we have $\sigma\left(
i\right)  <\sigma\left(  v\right)  $. In other words, $\sigma\left(  v\right)
>\sigma\left(  i\right)  $. Hence, $\left[  \sigma\left(  v\right)
>\sigma\left(  i\right)  \right]  =1$. Also, we don't have $\sigma\left(
i\right)  >\sigma\left(  v\right)  $ (since we have $\sigma\left(  i\right)
<\sigma\left(  v\right)  $). Hence, we have $\left[  \sigma\left(  i\right)
>\sigma\left(  v\right)  \right]  =0$. Moreover, we have $\sigma\left(
v\right)  >\sigma\left(  i\right)  >\sigma\left(  u\right)  $. Thus, $\left[
\sigma\left(  v\right)  >\sigma\left(  i\right)  >\sigma\left(  u\right)
\right]  =1$. Now, comparing%
\begin{align*}
&  \underbrace{\left[  \sigma\left(  u\right)  >\sigma\left(  i\right)
\right]  }_{=0}-\underbrace{\left[  \sigma\left(  v\right)  >\sigma\left(
i\right)  \right]  }_{=1}+\underbrace{\left[  \sigma\left(  i\right)
>\sigma\left(  v\right)  \right]  }_{=0}-\underbrace{\left[  \sigma\left(
i\right)  >\sigma\left(  u\right)  \right]  }_{=1}\\
&  =0-1+0-1=-2
\end{align*}
with%
\[
2\underbrace{\left[  \sigma\left(  u\right)  >\sigma\left(  i\right)
>\sigma\left(  v\right)  \right]  }_{=0}-2\underbrace{\left[  \sigma\left(
v\right)  >\sigma\left(  i\right)  >\sigma\left(  u\right)  \right]  }%
_{=1}=2\cdot0-2\cdot1=-2,
\]
we obtain%
\begin{align*}
&  \left[  \sigma\left(  u\right)  >\sigma\left(  i\right)  \right]  -\left[
\sigma\left(  v\right)  >\sigma\left(  i\right)  \right]  +\left[
\sigma\left(  i\right)  >\sigma\left(  v\right)  \right]  -\left[
\sigma\left(  i\right)  >\sigma\left(  u\right)  \right] \\
&  =2\left[  \sigma\left(  u\right)  >\sigma\left(  i\right)  >\sigma\left(
v\right)  \right]  -2\left[  \sigma\left(  v\right)  >\sigma\left(  i\right)
>\sigma\left(  u\right)  \right]  .
\end{align*}
Thus, Lemma \ref{lem.sol.perm.lisitij.uiv} is proven in Subcase 1.1.

Let us now consider Subcase 1.2. In this subcase, we have $\sigma\left(
i\right)  >\sigma\left(  v\right)  $. In other words, $\sigma\left(  v\right)
<\sigma\left(  i\right)  $. Hence, we don't have $\sigma\left(  v\right)
>\sigma\left(  i\right)  $. Thus, we have $\left[  \sigma\left(  v\right)
>\sigma\left(  i\right)  \right]  =0$. Also, $\left[  \sigma\left(  i\right)
>\sigma\left(  v\right)  \right]  =1$ (since $\sigma\left(  i\right)
>\sigma\left(  v\right)  $). Moreover, we don't have $\sigma\left(  v\right)
>\sigma\left(  i\right)  >\sigma\left(  u\right)  $ (since we don't have
$\sigma\left(  v\right)  >\sigma\left(  i\right)  $). Thus, $\left[
\sigma\left(  v\right)  >\sigma\left(  i\right)  >\sigma\left(  u\right)
\right]  =0$. Now, comparing%
\begin{align*}
&  \underbrace{\left[  \sigma\left(  u\right)  >\sigma\left(  i\right)
\right]  }_{=0}-\underbrace{\left[  \sigma\left(  v\right)  >\sigma\left(
i\right)  \right]  }_{=0}+\underbrace{\left[  \sigma\left(  i\right)
>\sigma\left(  v\right)  \right]  }_{=1}-\underbrace{\left[  \sigma\left(
i\right)  >\sigma\left(  u\right)  \right]  }_{=1}\\
&  =0-0+1-1=0
\end{align*}
with%
\[
2\underbrace{\left[  \sigma\left(  u\right)  >\sigma\left(  i\right)
>\sigma\left(  v\right)  \right]  }_{=0}-2\underbrace{\left[  \sigma\left(
v\right)  >\sigma\left(  i\right)  >\sigma\left(  u\right)  \right]  }%
_{=0}=2\cdot0-2\cdot0=0,
\]
we obtain%
\begin{align*}
&  \left[  \sigma\left(  u\right)  >\sigma\left(  i\right)  \right]  -\left[
\sigma\left(  v\right)  >\sigma\left(  i\right)  \right]  +\left[
\sigma\left(  i\right)  >\sigma\left(  v\right)  \right]  -\left[
\sigma\left(  i\right)  >\sigma\left(  u\right)  \right] \\
&  =2\left[  \sigma\left(  u\right)  >\sigma\left(  i\right)  >\sigma\left(
v\right)  \right]  -2\left[  \sigma\left(  v\right)  >\sigma\left(  i\right)
>\sigma\left(  u\right)  \right]  .
\end{align*}
Thus, Lemma \ref{lem.sol.perm.lisitij.uiv} is proven in Subcase 1.2.

We have now proven Lemma \ref{lem.sol.perm.lisitij.uiv} in each of the two
Subcases 1.1 and 1.2. Since these two Subcases cover the whole Case 1, we thus
conclude that Lemma \ref{lem.sol.perm.lisitij.uiv} holds in Case 1.

Let us now consider Case 2. In this case, we have $\sigma\left(  u\right)
>\sigma\left(  i\right)  $. In other words, $\sigma\left(  i\right)
<\sigma\left(  u\right)  $. Hence, we don't have $\sigma\left(  i\right)
>\sigma\left(  u\right)  $. Thus, we have $\left[  \sigma\left(  i\right)
>\sigma\left(  u\right)  \right]  =0$. Also, $\left[  \sigma\left(  u\right)
>\sigma\left(  i\right)  \right]  =1$ (since $\sigma\left(  u\right)
>\sigma\left(  i\right)  $). Also, we don't have $\sigma\left(  v\right)
>\sigma\left(  i\right)  >\sigma\left(  u\right)  $ (because we don't have
$\sigma\left(  i\right)  >\sigma\left(  u\right)  $). Thus, $\left[
\sigma\left(  v\right)  >\sigma\left(  i\right)  >\sigma\left(  u\right)
\right]  =0$.

We have $\sigma\left(  i\right)  \neq\sigma\left(  v\right)  $. In other
words, either $\sigma\left(  i\right)  <\sigma\left(  v\right)  $ or
$\sigma\left(  i\right)  >\sigma\left(  v\right)  $. Thus, we are in one of
the following two subcases:

\textit{Subcase 2.1:} We have $\sigma\left(  i\right)  <\sigma\left(
v\right)  $.

\textit{Subcase 2.2:} We have $\sigma\left(  i\right)  >\sigma\left(
v\right)  $.

Let us first consider Subcase 2.1. In this subcase, we have $\sigma\left(
i\right)  <\sigma\left(  v\right)  $. In other words, $\sigma\left(  v\right)
>\sigma\left(  i\right)  $. Hence, $\left[  \sigma\left(  v\right)
>\sigma\left(  i\right)  \right]  =1$. Also, we don't have $\sigma\left(
i\right)  >\sigma\left(  v\right)  $ (since we have $\sigma\left(  i\right)
<\sigma\left(  v\right)  $). Hence, we have $\left[  \sigma\left(  i\right)
>\sigma\left(  v\right)  \right]  =0$. Moreover, we don't have $\sigma\left(
u\right)  >\sigma\left(  i\right)  >\sigma\left(  v\right)  $ (since we don't
have $\sigma\left(  i\right)  >\sigma\left(  v\right)  $). Thus, $\left[
\sigma\left(  u\right)  >\sigma\left(  i\right)  >\sigma\left(  v\right)
\right]  =0$. Now, comparing%
\begin{align*}
&  \underbrace{\left[  \sigma\left(  u\right)  >\sigma\left(  i\right)
\right]  }_{=1}-\underbrace{\left[  \sigma\left(  v\right)  >\sigma\left(
i\right)  \right]  }_{=1}+\underbrace{\left[  \sigma\left(  i\right)
>\sigma\left(  v\right)  \right]  }_{=0}-\underbrace{\left[  \sigma\left(
i\right)  >\sigma\left(  u\right)  \right]  }_{=0}\\
&  =1-1+0-0=0
\end{align*}
with%
\[
2\underbrace{\left[  \sigma\left(  u\right)  >\sigma\left(  i\right)
>\sigma\left(  v\right)  \right]  }_{=0}-2\underbrace{\left[  \sigma\left(
v\right)  >\sigma\left(  i\right)  >\sigma\left(  u\right)  \right]  }%
_{=0}=2\cdot0-2\cdot0=0,
\]
we obtain%
\begin{align*}
&  \left[  \sigma\left(  u\right)  >\sigma\left(  i\right)  \right]  -\left[
\sigma\left(  v\right)  >\sigma\left(  i\right)  \right]  +\left[
\sigma\left(  i\right)  >\sigma\left(  v\right)  \right]  -\left[
\sigma\left(  i\right)  >\sigma\left(  u\right)  \right] \\
&  =2\left[  \sigma\left(  u\right)  >\sigma\left(  i\right)  >\sigma\left(
v\right)  \right]  -2\left[  \sigma\left(  v\right)  >\sigma\left(  i\right)
>\sigma\left(  u\right)  \right]  .
\end{align*}
Thus, Lemma \ref{lem.sol.perm.lisitij.uiv} is proven in Subcase 2.1.

Let us now consider Subcase 2.2. In this subcase, we have $\sigma\left(
i\right)  >\sigma\left(  v\right)  $. In other words, $\sigma\left(  v\right)
<\sigma\left(  i\right)  $. Hence, we don't have $\sigma\left(  v\right)
>\sigma\left(  i\right)  $. Thus, we have $\left[  \sigma\left(  v\right)
>\sigma\left(  i\right)  \right]  =0$. Also, $\left[  \sigma\left(  i\right)
>\sigma\left(  v\right)  \right]  =1$ (since $\sigma\left(  i\right)
>\sigma\left(  v\right)  $). Moreover, we have $\sigma\left(  u\right)
>\sigma\left(  i\right)  >\sigma\left(  v\right)  $. Thus, $\left[
\sigma\left(  u\right)  >\sigma\left(  i\right)  >\sigma\left(  v\right)
\right]  =1$. Now, comparing%
\begin{align*}
&  \underbrace{\left[  \sigma\left(  u\right)  >\sigma\left(  i\right)
\right]  }_{=1}-\underbrace{\left[  \sigma\left(  v\right)  >\sigma\left(
i\right)  \right]  }_{=0}+\underbrace{\left[  \sigma\left(  i\right)
>\sigma\left(  v\right)  \right]  }_{=1}-\underbrace{\left[  \sigma\left(
i\right)  >\sigma\left(  u\right)  \right]  }_{=0}\\
&  =1-0+1-0=2
\end{align*}
with%
\[
2\underbrace{\left[  \sigma\left(  u\right)  >\sigma\left(  i\right)
>\sigma\left(  v\right)  \right]  }_{=1}-2\underbrace{\left[  \sigma\left(
v\right)  >\sigma\left(  i\right)  >\sigma\left(  u\right)  \right]  }%
_{=0}=2\cdot1-2\cdot0=2,
\]
we obtain%
\begin{align*}
&  \left[  \sigma\left(  u\right)  >\sigma\left(  i\right)  \right]  -\left[
\sigma\left(  v\right)  >\sigma\left(  i\right)  \right]  +\left[
\sigma\left(  i\right)  >\sigma\left(  v\right)  \right]  -\left[
\sigma\left(  i\right)  >\sigma\left(  u\right)  \right] \\
&  =2\left[  \sigma\left(  u\right)  >\sigma\left(  i\right)  >\sigma\left(
v\right)  \right]  -2\left[  \sigma\left(  v\right)  >\sigma\left(  i\right)
>\sigma\left(  u\right)  \right]  .
\end{align*}
Thus, Lemma \ref{lem.sol.perm.lisitij.uiv} is proven in Subcase 2.2.

We have now proven Lemma \ref{lem.sol.perm.lisitij.uiv} in each of the two
Subcases 2.1 and 2.2. Since these two Subcases cover the whole Case 2, we thus
conclude that Lemma \ref{lem.sol.perm.lisitij.uiv} holds in Case 2.

We have now proven Lemma \ref{lem.sol.perm.lisitij.uiv} in each of the two
Cases 1 and 2. Since these two Cases cover all possibilities, we thus conclude
that Lemma \ref{lem.sol.perm.lisitij.uiv} always holds. This completes the
proof of Lemma \ref{lem.sol.perm.lisitij.uiv}.
\end{proof}
\end{verlong}

\begin{lemma}
\label{lem.sol.perm.lisitij.1}Let $n\in\mathbb{N}$. Let $u\in\left[  n\right]
$ and $v\in\left[  n\right]  $ be such that $u<v$. Let $\sigma\in S_{n}$ and
$\tau\in S_{n}$ be such that $\tau=\sigma\circ t_{u,v}$. Then:

\textbf{(a)} We have $\ell_{i}\left(  \sigma\right)  =\ell_{i}\left(
\tau\right)  $ for each $i\in\left\{  1,2,\ldots,u-1\right\}  $.

\textbf{(b)} We have%
\[
\ell_{u}\left(  \sigma\right)  =\ell_{v}\left(  \tau\right)  +\sum
_{i=u+1}^{v-1}\left[  \sigma\left(  u\right)  >\sigma\left(  i\right)
\right]  +\left[  \sigma\left(  u\right)  >\sigma\left(  v\right)  \right]  .
\]

\textbf{(c)} We have%
\[
\ell_{i}\left(  \sigma\right)  =\ell_{i}\left(  \tau\right)  +\left[
\sigma\left(  i\right)  >\sigma\left(  v\right)  \right]  -\left[
\sigma\left(  i\right)  >\sigma\left(  u\right)  \right]
\]
for each $i\in\left\{  u+1,u+2,\ldots,v-1\right\}  $.

\textbf{(d)} We have%
\[
\ell_{v}\left(  \sigma\right)  =\ell_{u}\left(  \tau\right)  -\sum
_{i=u+1}^{v-1}\left[  \sigma\left(  v\right)  >\sigma\left(  i\right)
\right]  -\left[  \sigma\left(  v\right)  >\sigma\left(  u\right)  \right]  .
\]

\textbf{(e)} We have $\ell_{i}\left(  \sigma\right)  =\ell_{i}\left(
\tau\right)  $ for each $i\in\left\{  v+1,v+2,\ldots,n\right\}  $.

\textbf{(f)} We have%
\begin{align*}
\ell\left(  \sigma\right)   &  =\ell\left(  \tau\right)  +\left[
\sigma\left(  u\right)  >\sigma\left(  v\right)  \right]  -\left[
\sigma\left(  v\right)  >\sigma\left(  u\right)  \right] \\
&  \ \ \ \ \ \ \ \ \ \ +\sum_{i=u+1}^{v-1}\left(  2\left[  \sigma\left(
u\right)  >\sigma\left(  i\right)  >\sigma\left(  v\right)  \right]  -2\left[
\sigma\left(  v\right)  >\sigma\left(  i\right)  >\sigma\left(  u\right)
\right]  \right)  .
\end{align*}

\end{lemma}

\begin{proof}
[Proof of Lemma \ref{lem.sol.perm.lisitij.1}.]We have $u\in\left[  n\right]
=\left\{  1,2,\ldots,n\right\}  $ (by the definition of $\left[  n\right]  $)
and $v\in\left\{  1,2,\ldots,n\right\}  $ (similarly). From $u\in\left\{
1,2,\ldots,n\right\}  $, we obtain $u\geq1$. From $v\in\left\{  1,2,\ldots
,n\right\}  $, we obtain $v\leq n$. Thus, $1\leq u<v\leq n$. Note that $u$ and
$v$ are distinct (since $u<v$). Also, from $u<v$, we obtain $u\leq v-1$ (since
$u$ and $v$ are integers), and thus $u+1\leq v$.

Recall that $t_{u,v}$ is the permutation in $S_{n}$ which swaps $u$ with $v$
while leaving all other elements of $\left\{  1,2,\ldots,n\right\}  $
unchanged (by the definition of $t_{u,v}$). In other words, $t_{u,v}$ is the
permutation in $S_{n}$ that satisfies $t_{u,v}\left(  u\right)  =v$ and
$t_{u,v}\left(  v\right)  =u$ and%
\begin{equation}
\left(  t_{u,v}\left(  i\right)  =i\ \ \ \ \ \ \ \ \ \ \text{for each }%
i\in\left\{  1,2,\ldots,n\right\}  \setminus\left\{  u,v\right\}  \right)  .
\label{pf.lem.sol.perm.lisitij.1.1}%
\end{equation}

Now,
\begin{align*}
\underbrace{\tau}_{=\sigma\circ t_{u,v}}\left(  u\right)   &  =\left(
\sigma\circ t_{u,v}\right)  \left(  u\right)  =\sigma\left(
\underbrace{t_{u,v}\left(  u\right)  }_{=v}\right)  =\sigma\left(  v\right)
\ \ \ \ \ \ \ \ \ \ \text{and}\\
\underbrace{\tau}_{=\sigma\circ t_{u,v}}\left(  v\right)   &  =\left(
\sigma\circ t_{u,v}\right)  \left(  v\right)  =\sigma\left(
\underbrace{t_{u,v}\left(  v\right)  }_{=u}\right)  =\sigma\left(  u\right)  .
\end{align*}
Moreover, each $i\in\left\{  1,2,\ldots,n\right\}  \setminus\left\{
u,v\right\}  $ satisfies%
\begin{equation}
\underbrace{\tau}_{=\sigma\circ t_{u,v}}\left(  i\right)  =\left(  \sigma\circ
t_{u,v}\right)  \left(  i\right)  =\sigma\left(  \underbrace{t_{u,v}\left(
i\right)  }_{\substack{=i\\\text{(by (\ref{pf.lem.sol.perm.lisitij.1.1}))}%
}}\right)  =\sigma\left(  i\right)  . \label{pf.lem.sol.perm.lisitij.1.ti=si}%
\end{equation}

From this, we obtain the following consequences:

\begin{itemize}
\item Each $i\in\left\{  1,2,\ldots,u-1\right\}  $ satisfies%
\begin{equation}
\tau\left(  i\right)  =\sigma\left(  i\right)  .
\label{pf.lem.sol.perm.lisitij.1.ti=si-left}%
\end{equation}

\begin{vershort}
[\textit{Proof of (\ref{pf.lem.sol.perm.lisitij.1.ti=si-left}):} Let
$i\in\left\{  1,2,\ldots,u-1\right\}  $. Thus, $i\leq u-1<u$, so that $i\neq
u$. Also, $i<u<v$, so that $i\neq v$. Also, $i\in\left\{  1,2,\ldots
,u-1\right\}  \subseteq\left\{  1,2,\ldots,n\right\}  $ (since $u-1\leq u\leq
n$). Combining $i\neq u$ with $i\neq v$, we conclude that $i\notin\left\{
u,v\right\}  $. Thus, $i\in\left\{  1,2,\ldots,n\right\}  \setminus\left\{
u,v\right\}  $ (since $i\in\left\{  1,2,\ldots,n\right\}  $). Hence,
(\ref{pf.lem.sol.perm.lisitij.1.ti=si}) yields $\tau\left(  i\right)
=\sigma\left(  i\right)  $. This proves
(\ref{pf.lem.sol.perm.lisitij.1.ti=si-left}).]
\end{vershort}

\begin{verlong}
[\textit{Proof of (\ref{pf.lem.sol.perm.lisitij.1.ti=si-left}):} Let
$i\in\left\{  1,2,\ldots,u-1\right\}  $. Thus, $i\leq u-1<u$, so that $i\neq
u$. Also, $i<u<v$, so that $i\neq v$. Combining $i\neq u$ with $i\neq v$, we
conclude that we have $\left(  \text{neither }i=u\text{ nor }i=v\right)  $.
Thus, $i\in\left\{  1,2,\ldots,n\right\}  \setminus\left\{  u,v\right\}
$\ \ \ \ \footnote{\textit{Proof.} If we had $i\in\left\{  u,v\right\}  $,
then we would have $\left(  \text{either }i=u\text{ or }i=v\right)  $, which
would contradict the fact that $\left(  \text{neither }i=u\text{ nor
}i=v\right)  $. Hence, we cannot have $i\in\left\{  u,v\right\}  $. In other
words, we have $i\notin\left\{  u,v\right\}  $. Also, $i\in\left\{
1,2,\ldots,u-1\right\}  \subseteq\left\{  1,2,\ldots,n\right\}  $ (since
$u-1<u\leq n$ (because $u\in\left\{  1,2,\ldots,n\right\}  $)). Combining this
with $i\notin\left\{  u,v\right\}  $, we obtain $i\in\left\{  1,2,\ldots
,n\right\}  \setminus\left\{  u,v\right\}  $. Qed.}. Hence,
(\ref{pf.lem.sol.perm.lisitij.1.ti=si}) yields $\tau\left(  i\right)
=\sigma\left(  i\right)  $. This proves
(\ref{pf.lem.sol.perm.lisitij.1.ti=si-left}).]
\end{verlong}

\item Each $i\in\left\{  u+1,u+2,\ldots,v-1\right\}  $ satisfies%
\begin{equation}
\tau\left(  i\right)  =\sigma\left(  i\right)  .
\label{pf.lem.sol.perm.lisitij.1.ti=si-mid}%
\end{equation}

\begin{vershort}
[\textit{Proof of (\ref{pf.lem.sol.perm.lisitij.1.ti=si-mid}):} Let
$i\in\left\{  u+1,u+2,\ldots,v-1\right\}  $. Thus, $u+1\leq i\leq v-1$. We
have $u+1\leq i$, thus $i\geq u+1>u$ and therefore $i\neq u$. Also, $i\leq
v-1<v$ and thus $i\neq v$. Moreover, $i>u\geq1$ and thus $i\geq1$. Also,
$i<v\leq n$ and thus $i\leq n$. Thus, $i\in\left\{  1,2,\ldots,n\right\}  $
(since $i\geq1$). Combining $i\neq u$ with $i\neq v$, we conclude that
$i\notin\left\{  u,v\right\}  $. Thus, $i\in\left\{  1,2,\ldots,n\right\}
\setminus\left\{  u,v\right\}  $ (since $i\in\left\{  1,2,\ldots,n\right\}
$). Hence, (\ref{pf.lem.sol.perm.lisitij.1.ti=si}) yields $\tau\left(
i\right)  =\sigma\left(  i\right)  $. This proves
(\ref{pf.lem.sol.perm.lisitij.1.ti=si-mid}).]
\end{vershort}

\begin{verlong}
[\textit{Proof of (\ref{pf.lem.sol.perm.lisitij.1.ti=si-mid}):} Let
$i\in\left\{  u+1,u+2,\ldots,v-1\right\}  $. Thus, $u+1\leq i\leq v-1$. We
have $u+1\leq i$, thus $i\geq u+1>u$ and therefore $i\neq u$. Also, $i\leq
v-1<v$ and thus $i\neq v$. Moreover, $i>u\geq1$ and thus $i\geq1$. Also,
$i<v\leq n$ and thus $i\leq n$. Combining $i\neq u$ with $i\neq v$, we
conclude that we have $\left(  \text{neither }i=u\text{ nor }i=v\right)  $.
Thus, $i\in\left\{  1,2,\ldots,n\right\}  \setminus\left\{  u,v\right\}
$\ \ \ \ \footnote{\textit{Proof.} If we had $i\in\left\{  u,v\right\}  $,
then we would have $\left(  \text{either }i=u\text{ or }i=v\right)  $, which
would contradict the fact that $\left(  \text{neither }i=u\text{ nor
}i=v\right)  $. Hence, we cannot have $i\in\left\{  u,v\right\}  $. In other
words, we have $i\notin\left\{  u,v\right\}  $. Also, $i\in\left\{
1,2,\ldots,n\right\}  $ (since $i\geq1$ and $i\leq n$). Combining this with
$i\notin\left\{  u,v\right\}  $, we obtain $i\in\left\{  1,2,\ldots,n\right\}
\setminus\left\{  u,v\right\}  $. Qed.}. Hence,
(\ref{pf.lem.sol.perm.lisitij.1.ti=si}) yields $\tau\left(  i\right)
=\sigma\left(  i\right)  $. This proves
(\ref{pf.lem.sol.perm.lisitij.1.ti=si-mid}).]
\end{verlong}

\item Each $j\in\left\{  v+1,v+2,\ldots,n\right\}  $ satisfies%
\begin{equation}
\tau\left(  j\right)  =\sigma\left(  j\right)  .
\label{pf.lem.sol.perm.lisitij.1.ti=si-right}%
\end{equation}

\begin{vershort}
[\textit{Proof of (\ref{pf.lem.sol.perm.lisitij.1.ti=si-right}):} Let
$j\in\left\{  v+1,v+2,\ldots,n\right\}  $. Thus, $j\geq v+1>v$, so that $j\neq
v$. Also, $j>v>u$ (since $u<v$), so that $j\neq u$. Also, $j\in\left\{
v+1,v+2,\ldots,n\right\}  \subseteq\left\{  1,2,\ldots,n\right\}  $ (since
$v+1\geq v\geq1$). Combining $j\neq u$ with $j\neq v$, we conclude that
$j\notin\left\{  u,v\right\}  $. Thus, $j\in\left\{  1,2,\ldots,n\right\}
\setminus\left\{  u,v\right\}  $ (since $j\in\left\{  1,2,\ldots,n\right\}
$). Hence, (\ref{pf.lem.sol.perm.lisitij.1.ti=si}) (applied to $i=j$) yields
$\tau\left(  j\right)  =\sigma\left(  j\right)  $. This proves
(\ref{pf.lem.sol.perm.lisitij.1.ti=si-right}).]
\end{vershort}

\begin{verlong}
[\textit{Proof of (\ref{pf.lem.sol.perm.lisitij.1.ti=si-right}):} Let
$j\in\left\{  v+1,v+2,\ldots,n\right\}  $. Thus, $j\geq v+1>v$, so that $j\neq
v$. Also, $j>v>u$ (since $u<v$), so that $j\neq u$. Combining $j\neq u$ with
$j\neq v$, we conclude that we have $\left(  \text{neither }j=u\text{ nor
}j=v\right)  $. Thus, $j\in\left\{  1,2,\ldots,n\right\}  \setminus\left\{
u,v\right\}  $\ \ \ \ \footnote{\textit{Proof.} If we had $j\in\left\{
u,v\right\}  $, then we would have $\left(  \text{either }j=u\text{ or
}j=v\right)  $, which would contradict the fact that $\left(  \text{neither
}j=u\text{ nor }j=v\right)  $. Hence, we cannot have $j\in\left\{
u,v\right\}  $. In other words, we have $j\notin\left\{  u,v\right\}  $. Also,
$j\in\left\{  v+1,v+2,\ldots,n\right\}  \subseteq\left\{  1,2,\ldots
,n\right\}  $ (since $v+1\geq v\geq1$ (because $v\in\left\{  1,2,\ldots
,n\right\}  $)). Combining this with $j\notin\left\{  u,v\right\}  $, we
obtain $j\in\left\{  1,2,\ldots,n\right\}  \setminus\left\{  u,v\right\}  $.
Qed.}. Hence, (\ref{pf.lem.sol.perm.lisitij.1.ti=si}) (applied to $i=j$)
yields $\tau\left(  j\right)  =\sigma\left(  j\right)  $. This proves
(\ref{pf.lem.sol.perm.lisitij.1.ti=si-right}).]
\end{verlong}
\end{itemize}

\textbf{(a)} Let $i\in\left\{  1,2,\ldots,u-1\right\}  $. Thus,
(\ref{pf.lem.sol.perm.lisitij.1.ti=si-left}) yields $\tau\left(  i\right)
=\sigma\left(  i\right)  $.

Also, each $j\in\left[  i\right]  $ yields%
\begin{equation}
\tau\left(  j\right)  =\sigma\left(  j\right)  .
\label{pf.lem.sol.perm.lisitij.1.a.tj=sj}%
\end{equation}

[\textit{Proof of (\ref{pf.lem.sol.perm.lisitij.1.a.tj=sj}):} Let $j\in\left[
i\right]  $. Thus, $j\in\left[  i\right]  =\left\{  1,2,\ldots,i\right\}  $
(by the definition of $\left[  i\right]  $). Hence, $j\geq1$ and $j\leq i\leq
u-1$ (since $i\in\left\{  1,2,\ldots,u-1\right\}  $), so that $j\in\left\{
1,2,\ldots,u-1\right\}  $ (because $j\geq1$). Thus,
(\ref{pf.lem.sol.perm.lisitij.1.ti=si-left}) (applied to $j$ instead of $i$)
yields $\tau\left(  j\right)  =\sigma\left(  j\right)  $. This proves
(\ref{pf.lem.sol.perm.lisitij.1.a.tj=sj}).]

Now,%
\[
\tau\left(  \left[  i\right]  \right)  =\left\{  \underbrace{\tau\left(
j\right)  }_{\substack{=\sigma\left(  j\right)  \\\text{(by
(\ref{pf.lem.sol.perm.lisitij.1.a.tj=sj}))}}}\ \mid\ j\in\left[  i\right]
\right\}  =\left\{  \sigma\left(  j\right)  \ \mid\ j\in\left[  i\right]
\right\}  =\sigma\left(  \left[  i\right]  \right)  .
\]

Lemma \ref{lem.perm.lexico1.lis} \textbf{(a)} yields $\ell_{i}\left(
\sigma\right)  =\left\vert \left[  \sigma\left(  i\right)  -1\right]
\setminus\sigma\left(  \left[  i\right]  \right)  \right\vert $. The same
argument (applied to $\tau$ instead of $\sigma$) yields $\ell_{i}\left(
\tau\right)  =\left\vert \left[  \tau\left(  i\right)  -1\right]
\setminus\tau\left(  \left[  i\right]  \right)  \right\vert $. Hence,%
\[
\ell_{i}\left(  \tau\right)  =\left\vert \left[  \underbrace{\tau\left(
i\right)  }_{=\sigma\left(  i\right)  }-1\right]  \setminus\underbrace{\tau
\left(  \left[  i\right]  \right)  }_{=\sigma\left(  \left[  i\right]
\right)  }\right\vert =\left\vert \left[  \sigma\left(  i\right)  -1\right]
\setminus\sigma\left(  \left[  i\right]  \right)  \right\vert =\ell_{i}\left(
\sigma\right)
\]
(since $\ell_{i}\left(  \sigma\right)  =\left\vert \left[  \sigma\left(
i\right)  -1\right]  \setminus\sigma\left(  \left[  i\right]  \right)
\right\vert $). In other words, $\ell_{i}\left(  \sigma\right)  =\ell
_{i}\left(  \tau\right)  $. This proves Lemma \ref{lem.sol.perm.lisitij.1}
\textbf{(a)}.

\textbf{(b)} Lemma \ref{lem.sol.perm.lisitij.li=} (applied to $i=u$) yields%
\begin{align}
\ell_{u}\left(  \sigma\right)   &  =\sum_{j=u+1}^{n}\left[  \sigma\left(
u\right)  >\sigma\left(  j\right)  \right]  =\sum_{j=u+1}^{v-1}\left[
\sigma\left(  u\right)  >\sigma\left(  j\right)  \right]  +\underbrace{\sum
_{j=v}^{n}\left[  \sigma\left(  u\right)  >\sigma\left(  j\right)  \right]
}_{\substack{=\left[  \sigma\left(  u\right)  >\sigma\left(  v\right)
\right]  +\sum_{j=v+1}^{n}\left[  \sigma\left(  u\right)  >\sigma\left(
j\right)  \right]  \\\text{(here, we have split off the addend for
}j=v\\\text{from the sum, since }v\leq n\text{)}}}\nonumber\\
&  \ \ \ \ \ \ \ \ \ \ \left(
\begin{array}
[c]{c}%
\text{here, we have split the sum at }j=v\text{,}\\
\text{because }u+1\leq v\leq n
\end{array}
\right) \nonumber\\
&  =\sum_{j=u+1}^{v-1}\left[  \sigma\left(  u\right)  >\sigma\left(  j\right)
\right]  +\left[  \sigma\left(  u\right)  >\sigma\left(  v\right)  \right]
+\sum_{j=v+1}^{n}\left[  \sigma\left(  u\right)  >\sigma\left(  j\right)
\right]  . \label{pf.lem.sol.perm.lisitij.1.b.1}%
\end{align}
But Lemma \ref{lem.sol.perm.lisitij.li=} (applied to $v$ and $\tau$ instead of
$i$ and $\sigma$) yields%
\begin{equation}
\ell_{v}\left(  \tau\right)  =\sum_{j=v+1}^{n}\left[  \underbrace{\tau\left(
v\right)  }_{=\sigma\left(  u\right)  }>\underbrace{\tau\left(  j\right)
}_{\substack{=\sigma\left(  j\right)  \\\text{(by
(\ref{pf.lem.sol.perm.lisitij.1.ti=si-right}))}}}\right]  =\sum_{j=v+1}%
^{n}\left[  \sigma\left(  u\right)  >\sigma\left(  j\right)  \right]  .
\label{pf.lem.sol.perm.lisitij.1.b.2}%
\end{equation}
Hence, (\ref{pf.lem.sol.perm.lisitij.1.b.1}) becomes%
\begin{align*}
\ell_{u}\left(  \sigma\right)   &  =\underbrace{\sum_{j=u+1}^{v-1}\left[
\sigma\left(  u\right)  >\sigma\left(  j\right)  \right]  }_{\substack{=\sum
_{i=u+1}^{v-1}\left[  \sigma\left(  u\right)  >\sigma\left(  i\right)
\right]  \\\text{(here, we have renamed}\\\text{the summation index }j\text{
as }i\text{)}}}+\left[  \sigma\left(  u\right)  >\sigma\left(  v\right)
\right]  +\underbrace{\sum_{j=v+1}^{n}\left[  \sigma\left(  u\right)
>\sigma\left(  j\right)  \right]  }_{\substack{=\ell_{v}\left(  \tau\right)
\\\text{(by (\ref{pf.lem.sol.perm.lisitij.1.b.2}))}}}\\
&  =\sum_{i=u+1}^{v-1}\left[  \sigma\left(  u\right)  >\sigma\left(  i\right)
\right]  +\left[  \sigma\left(  u\right)  >\sigma\left(  v\right)  \right]
+\ell_{v}\left(  \tau\right) \\
&  =\ell_{v}\left(  \tau\right)  +\sum_{i=u+1}^{v-1}\left[  \sigma\left(
u\right)  >\sigma\left(  i\right)  \right]  +\left[  \sigma\left(  u\right)
>\sigma\left(  v\right)  \right]  .
\end{align*}
This proves Lemma \ref{lem.sol.perm.lisitij.1} \textbf{(b)}.

\textbf{(c)} Let $i\in\left\{  u+1,u+2,\ldots,v-1\right\}  $. Thus, $u+1\leq
i\leq v-1$. Hence, $v\in\left\{  i+1,i+2,\ldots,n\right\}  $%
\ \ \ \ \footnote{\textit{Proof.} We have $i\leq v-1$, so that $v\geq i+1$.
Combining this with $v\leq n$, we obtain $v\in\left\{  i+1,i+2,\ldots
,n\right\}  $ (since $v$ is an integer).}.

Furthermore, each $j\in\left\{  i+1,i+2,\ldots,n\right\}  $ satisfying $j\neq
v$ must satisfy%
\begin{equation}
\tau\left(  j\right)  =\sigma\left(  j\right)  .
\label{pf.lem.sol.perm.lisitij.1.c.tj=sj}%
\end{equation}

\begin{vershort}
[\textit{Proof of (\ref{pf.lem.sol.perm.lisitij.1.c.tj=sj}):} Let
$j\in\left\{  i+1,i+2,\ldots,n\right\}  $ be such that $j\neq v$.

From $j\in\left\{  i+1,i+2,\ldots,n\right\}  $, we obtain $j\geq i+1>i>u$
(since $u<u+1\leq i$) and therefore $j\neq u$. Also, $j\in\left\{
i+1,i+2,\ldots,n\right\}  \subseteq\left\{  1,2,\ldots,n\right\}  $ (since
$i+1\geq i>u\geq1$). Combining $j\neq u$ with $j\neq v$, we conclude that
$j\notin\left\{  u,v\right\}  $. Thus, $j\in\left\{  1,2,\ldots,n\right\}
\setminus\left\{  u,v\right\}  $ (since $j\in\left\{  1,2,\ldots,n\right\}
$). Hence, (\ref{pf.lem.sol.perm.lisitij.1.ti=si}) (applied to $j$ instead of
$i$) yields $\tau\left(  j\right)  =\sigma\left(  j\right)  $. This proves
(\ref{pf.lem.sol.perm.lisitij.1.c.tj=sj}).]
\end{vershort}

\begin{verlong}
[\textit{Proof of (\ref{pf.lem.sol.perm.lisitij.1.c.tj=sj}):} Let
$j\in\left\{  i+1,i+2,\ldots,n\right\}  $ be such that $j\neq v$. We must
prove (\ref{pf.lem.sol.perm.lisitij.1.c.tj=sj}).

We have $j\in\left\{  i+1,i+2,\ldots,n\right\}  $, thus $j\geq i+1>i>u$ (since
$u<u+1\leq i$) and therefore $j\neq u$. Combining $j\neq u$ with $j\neq v$, we
conclude that we have $\left(  \text{neither }j=u\text{ nor }j=v\right)  $.
Thus, $j\in\left\{  1,2,\ldots,n\right\}  \setminus\left\{  u,v\right\}
$\ \ \ \ \footnote{\textit{Proof.} If we had $j\in\left\{  u,v\right\}  $,
then we would have $\left(  \text{either }j=u\text{ or }j=v\right)  $, which
would contradict the fact that $\left(  \text{neither }j=u\text{ nor
}j=v\right)  $. Hence, we cannot have $j\in\left\{  u,v\right\}  $. In other
words, we have $j\notin\left\{  u,v\right\}  $. Also, $j\in\left\{
i+1,i+2,\ldots,n\right\}  \subseteq\left\{  1,2,\ldots,n\right\}  $ (since
$i+1\geq i>u\geq1$). Combining this with $j\notin\left\{  u,v\right\}  $, we
obtain $j\in\left\{  1,2,\ldots,n\right\}  \setminus\left\{  u,v\right\}  $.
Qed.}. Hence, (\ref{pf.lem.sol.perm.lisitij.1.ti=si}) (applied to $j$ instead
of $i$) yields $\tau\left(  j\right)  =\sigma\left(  j\right)  $. This proves
(\ref{pf.lem.sol.perm.lisitij.1.c.tj=sj}).]
\end{verlong}

Lemma \ref{lem.sol.perm.lisitij.li=} yields%
\begin{align}
\ell_{i}\left(  \sigma\right)   &  =\underbrace{\sum_{j=i+1}^{n}}_{=\sum
_{j\in\left\{  i+1,i+2,\ldots,n\right\}  }}\left[  \sigma\left(  i\right)
>\sigma\left(  j\right)  \right]  =\sum_{j\in\left\{  i+1,i+2,\ldots
,n\right\}  }\left[  \sigma\left(  i\right)  >\sigma\left(  j\right)  \right]
\nonumber\\
&  =\left[  \sigma\left(  i\right)  >\sigma\left(  v\right)  \right]
+\sum_{\substack{j\in\left\{  i+1,i+2,\ldots,n\right\}  ;\\j\neq v}}\left[
\sigma\left(  i\right)  >\sigma\left(  j\right)  \right]
\label{pf.lem.sol.perm.lisitij.1.c.sig}%
\end{align}
(here, we have split off the addend for $j=v$ from the sum, since we have
$v\in\left\{  i+1,i+2,\ldots,n\right\}  $). The same argument (but applied to
$\tau$ instead of $\sigma$) yields%
\[
\ell_{i}\left(  \tau\right)  =\left[  \tau\left(  i\right)  >\tau\left(
v\right)  \right]  +\sum_{\substack{j\in\left\{  i+1,i+2,\ldots,n\right\}
;\\j\neq v}}\left[  \tau\left(  i\right)  >\tau\left(  j\right)  \right]  .
\]
Thus,%
\begin{align*}
\ell_{i}\left(  \tau\right)   &  =\left[  \underbrace{\tau\left(  i\right)
}_{\substack{=\sigma\left(  i\right)  \\\text{(by
(\ref{pf.lem.sol.perm.lisitij.1.ti=si-mid}))}}}>\underbrace{\tau\left(
v\right)  }_{=\sigma\left(  u\right)  }\right]  +\sum_{\substack{j\in\left\{
i+1,i+2,\ldots,n\right\}  ;\\j\neq v}}\left[  \underbrace{\tau\left(
i\right)  }_{\substack{=\sigma\left(  i\right)  \\\text{(by
(\ref{pf.lem.sol.perm.lisitij.1.ti=si-mid}))}}}>\underbrace{\tau\left(
j\right)  }_{\substack{=\sigma\left(  j\right)  \\\text{(by
(\ref{pf.lem.sol.perm.lisitij.1.c.tj=sj}))}}}\right] \\
&  =\left[  \sigma\left(  i\right)  >\sigma\left(  u\right)  \right]
+\sum_{\substack{j\in\left\{  i+1,i+2,\ldots,n\right\}  ;\\j\neq v}}\left[
\sigma\left(  i\right)  >\sigma\left(  j\right)  \right]  .
\end{align*}
Subtracting this equality from (\ref{pf.lem.sol.perm.lisitij.1.c.sig}), we
obtain%
\begin{align*}
\ell_{i}\left(  \sigma\right)  -\ell_{i}\left(  \tau\right)   &  =\left(
\left[  \sigma\left(  i\right)  >\sigma\left(  v\right)  \right]
+\sum_{\substack{j\in\left\{  i+1,i+2,\ldots,n\right\}  ;\\j\neq v}}\left[
\sigma\left(  i\right)  >\sigma\left(  j\right)  \right]  \right) \\
&  \ \ \ \ \ \ \ \ \ \ -\left(  \left[  \sigma\left(  i\right)  >\sigma\left(
u\right)  \right]  +\sum_{\substack{j\in\left\{  i+1,i+2,\ldots,n\right\}
;\\j\neq v}}\left[  \sigma\left(  i\right)  >\sigma\left(  j\right)  \right]
\right) \\
&  =\left[  \sigma\left(  i\right)  >\sigma\left(  v\right)  \right]  -\left[
\sigma\left(  i\right)  >\sigma\left(  u\right)  \right]  .
\end{align*}
Adding $\ell_{i}\left(  \tau\right)  $ to both sides of this equality, we find%
\[
\ell_{i}\left(  \sigma\right)  =\ell_{i}\left(  \tau\right)  +\left[
\sigma\left(  i\right)  >\sigma\left(  v\right)  \right]  -\left[
\sigma\left(  i\right)  >\sigma\left(  u\right)  \right]  .
\]
This proves Lemma \ref{lem.sol.perm.lisitij.1} \textbf{(c)}.

\textbf{(d)} Lemma \ref{lem.sol.transpose.code.tij1} \textbf{(d)} (applied to
$u$ and $v$ instead of $i$ and $j$) yields $t_{u,v}\circ t_{u,v}%
=\operatorname*{id}$. Now,%
\[
\underbrace{\tau}_{=\sigma\circ t_{u,v}}\circ t_{u,v}=\sigma\circ
\underbrace{t_{u,v}\circ t_{u,v}}_{=\operatorname*{id}}=\sigma,
\]
so that $\sigma=\tau\circ t_{u,v}$. Hence, we can apply Lemma
\ref{lem.sol.perm.lisitij.1} \textbf{(b)} to $\sigma$ and $\tau$ instead of
$\tau$ and $\sigma$. We thus find%
\begin{align*}
\ell_{u}\left(  \tau\right)   &  =\ell_{v}\left(  \sigma\right)  +\sum
_{i=u+1}^{v-1}\left[  \underbrace{\tau\left(  u\right)  }_{=\sigma\left(
v\right)  }>\underbrace{\tau\left(  i\right)  }_{\substack{=\sigma\left(
i\right)  \\\text{(by (\ref{pf.lem.sol.perm.lisitij.1.ti=si-mid}))}}}\right]
+\left[  \underbrace{\tau\left(  u\right)  }_{=\sigma\left(  v\right)
}>\underbrace{\tau\left(  v\right)  }_{=\sigma\left(  u\right)  }\right] \\
&  =\ell_{v}\left(  \sigma\right)  +\sum_{i=u+1}^{v-1}\left[  \sigma\left(
v\right)  >\sigma\left(  i\right)  \right]  +\left[  \sigma\left(  v\right)
>\sigma\left(  u\right)  \right]  .
\end{align*}
Solving this equation for $\ell_{v}\left(  \sigma\right)  $, we obtain%
\[
\ell_{v}\left(  \sigma\right)  =\ell_{u}\left(  \tau\right)  -\sum
_{i=u+1}^{v-1}\left[  \sigma\left(  v\right)  >\sigma\left(  i\right)
\right]  -\left[  \sigma\left(  v\right)  >\sigma\left(  u\right)  \right]  .
\]
This proves Lemma \ref{lem.sol.perm.lisitij.1} \textbf{(d)}.

\textbf{(e)} Let $i\in\left\{  v+1,v+2,\ldots,n\right\}  $. Then,
(\ref{pf.lem.sol.perm.lisitij.1.ti=si-right}) (applied to $j=i$) yields
$\tau\left(  i\right)  =\sigma\left(  i\right)  $. Also, each $j\in\left\{
i+1,i+2,\ldots,n\right\}  $ satisfies%
\begin{equation}
\tau\left(  j\right)  =\sigma\left(  j\right)  .
\label{pf.lem.sol.perm.lisitij.1.e.tj=sj}%
\end{equation}

[\textit{Proof of (\ref{pf.lem.sol.perm.lisitij.1.e.tj=sj}):} Let
$j\in\left\{  i+1,i+2,\ldots,n\right\}  $. Thus, $j\geq i+1>i\geq v+1$ (since
$i\in\left\{  v+1,v+2,\ldots,n\right\}  $), so that $j\geq v+1$. Also, $j\leq
n$ (since $j\in\left\{  i+1,i+2,\ldots,n\right\}  $). Combining $j\geq v+1$
with $j\leq n$, we obtain $j\in\left\{  v+1,v+2,\ldots,n\right\}  $. Thus,
(\ref{pf.lem.sol.perm.lisitij.1.ti=si-right}) yields $\tau\left(  j\right)
=\sigma\left(  j\right)  $. This proves
(\ref{pf.lem.sol.perm.lisitij.1.e.tj=sj}).]

Lemma \ref{lem.sol.perm.lisitij.li=} yields%
\begin{equation}
\ell_{i}\left(  \sigma\right)  =\sum_{j=i+1}^{n}\left[  \sigma\left(
i\right)  >\sigma\left(  j\right)  \right]  .\nonumber
\end{equation}
Lemma \ref{lem.sol.perm.lisitij.li=} (applied to $\tau$ instead of $\sigma$)
yields%
\begin{equation}
\ell_{i}\left(  \tau\right)  =\sum_{j=i+1}^{n}\left[  \underbrace{\tau\left(
i\right)  }_{=\sigma\left(  i\right)  }>\underbrace{\tau\left(  j\right)
}_{\substack{=\sigma\left(  j\right)  \\\text{(by
(\ref{pf.lem.sol.perm.lisitij.1.e.tj=sj}))}}}\right]  =\sum_{j=i+1}^{n}\left[
\sigma\left(  i\right)  >\sigma\left(  j\right)  \right]  .\nonumber
\end{equation}
Comparing these two equalities, we obtain $\ell_{i}\left(  \sigma\right)
=\ell_{i}\left(  \tau\right)  $. This proves Lemma
\ref{lem.sol.perm.lisitij.1} \textbf{(e)}.

\textbf{(f)} Proposition \ref{prop.perm.lehmer.l} yields%
\begin{align}
\ell\left(  \sigma\right)   &  =\ell_{1}\left(  \sigma\right)  +\ell
_{2}\left(  \sigma\right)  +\cdots+\ell_{n}\left(  \sigma\right)  =\sum
_{i=1}^{n}\ell_{i}\left(  \sigma\right) \nonumber\\
&  =\underbrace{\sum_{i=1}^{v-1}\ell_{i}\left(  \sigma\right)  }%
_{\substack{=\sum_{i=1}^{u-1}\ell_{i}\left(  \sigma\right)  +\sum_{i=u}%
^{v-1}\ell_{i}\left(  \sigma\right)  \\\text{(here, we have split the sum at
}i=u\text{,}\\\text{since }1\leq u\leq v-1\text{)}}}+\underbrace{\sum
_{i=v}^{n}\ell_{i}\left(  \sigma\right)  }_{\substack{=\ell_{v}\left(
\sigma\right)  +\sum_{i=v+1}^{n}\ell_{i}\left(  \sigma\right)  \\\text{(here,
we have split off the addend for }i=v\text{ from the sum,}\\\text{since }v\leq
n\text{)}}}\nonumber\\
&  \ \ \ \ \ \ \ \ \ \ \left(  \text{here, we have split the sum at
}i=v\text{, because }1\leq v\leq n\right) \nonumber\\
&  =\sum_{i=1}^{u-1}\ell_{i}\left(  \sigma\right)  +\underbrace{\sum
_{i=u}^{v-1}\ell_{i}\left(  \sigma\right)  }_{\substack{=\ell_{u}\left(
\sigma\right)  +\sum_{i=u+1}^{v-1}\ell_{i}\left(  \sigma\right)
\\\text{(here, we have split off the addend for }i=u\text{ from the
sum,}\\\text{since }u\leq v-1\text{)}}}+\ell_{v}\left(  \sigma\right)
+\sum_{i=v+1}^{n}\ell_{i}\left(  \sigma\right) \nonumber\\
&  =\sum_{i=1}^{u-1}\ell_{i}\left(  \sigma\right)  +\ell_{u}\left(
\sigma\right)  +\sum_{i=u+1}^{v-1}\ell_{i}\left(  \sigma\right)  +\ell
_{v}\left(  \sigma\right)  +\sum_{i=v+1}^{n}\ell_{i}\left(  \sigma\right)  .
\label{pf.lem.sol.perm.lisitij.1.f.lsig=}%
\end{align}
The same argument (applied to $\tau$ instead of $\sigma$) yields%
\begin{equation}
\ell\left(  \tau\right)  =\sum_{i=1}^{u-1}\ell_{i}\left(  \tau\right)
+\ell_{u}\left(  \tau\right)  +\sum_{i=u+1}^{v-1}\ell_{i}\left(  \tau\right)
+\ell_{v}\left(  \tau\right)  +\sum_{i=v+1}^{n}\ell_{i}\left(  \tau\right)  .
\label{pf.lem.sol.perm.lisitij.1.f.ltau=}%
\end{equation}
But (\ref{pf.lem.sol.perm.lisitij.1.f.lsig=}) becomes%
\begin{align*}
\ell\left(  \sigma\right)   &  =\sum_{i=1}^{u-1}\underbrace{\ell_{i}\left(
\sigma\right)  }_{\substack{=\ell_{i}\left(  \tau\right)  \\\text{(by Lemma
\ref{lem.sol.perm.lisitij.1} \textbf{(a)})}}}+\underbrace{\ell_{u}\left(
\sigma\right)  }_{\substack{=\ell_{v}\left(  \tau\right)  +\sum_{i=u+1}%
^{v-1}\left[  \sigma\left(  u\right)  >\sigma\left(  i\right)  \right]
+\left[  \sigma\left(  u\right)  >\sigma\left(  v\right)  \right]  \\\text{(by
Lemma \ref{lem.sol.perm.lisitij.1} \textbf{(b)})}}}\\
&  \ \ \ \ \ \ \ \ \ \ +\sum_{i=u+1}^{v-1}\underbrace{\ell_{i}\left(
\sigma\right)  }_{\substack{=\ell_{i}\left(  \tau\right)  +\left[
\sigma\left(  i\right)  >\sigma\left(  v\right)  \right]  -\left[
\sigma\left(  i\right)  >\sigma\left(  u\right)  \right]  \\\text{(by Lemma
\ref{lem.sol.perm.lisitij.1} \textbf{(c)})}}}\\
&  \ \ \ \ \ \ \ \ \ \ +\underbrace{\ell_{v}\left(  \sigma\right)
}_{\substack{=\ell_{u}\left(  \tau\right)  -\sum_{i=u+1}^{v-1}\left[
\sigma\left(  v\right)  >\sigma\left(  i\right)  \right]  -\left[
\sigma\left(  v\right)  >\sigma\left(  u\right)  \right]  \\\text{(by Lemma
\ref{lem.sol.perm.lisitij.1} \textbf{(d)})}}}+\sum_{i=v+1}^{n}\underbrace{\ell
_{i}\left(  \sigma\right)  }_{\substack{=\ell_{i}\left(  \tau\right)
\\\text{(by Lemma \ref{lem.sol.perm.lisitij.1} \textbf{(e)})}}}\\
&  =\sum_{i=1}^{u-1}\ell_{i}\left(  \tau\right)  +\ell_{v}\left(  \tau\right)
+\sum_{i=u+1}^{v-1}\left[  \sigma\left(  u\right)  >\sigma\left(  i\right)
\right]  +\left[  \sigma\left(  u\right)  >\sigma\left(  v\right)  \right] \\
&  \ \ \ \ \ \ \ \ \ \ +\underbrace{\sum_{i=u+1}^{v-1}\left(  \ell_{i}\left(
\tau\right)  +\left[  \sigma\left(  i\right)  >\sigma\left(  v\right)
\right]  -\left[  \sigma\left(  i\right)  >\sigma\left(  u\right)  \right]
\right)  }_{=\sum_{i=u+1}^{v-1}\ell_{i}\left(  \tau\right)  +\sum
_{i=u+1}^{v-1}\left(  \left[  \sigma\left(  i\right)  >\sigma\left(  v\right)
\right]  -\left[  \sigma\left(  i\right)  >\sigma\left(  u\right)  \right]
\right)  }\\
&  \ \ \ \ \ \ \ \ \ \ +\ell_{u}\left(  \tau\right)  -\sum_{i=u+1}%
^{v-1}\left[  \sigma\left(  v\right)  >\sigma\left(  i\right)  \right]
-\left[  \sigma\left(  v\right)  >\sigma\left(  u\right)  \right]
+\sum_{i=v+1}^{n}\ell_{i}\left(  \tau\right) \\
&  =\sum_{i=1}^{u-1}\ell_{i}\left(  \tau\right)  +\ell_{v}\left(  \tau\right)
+\sum_{i=u+1}^{v-1}\left[  \sigma\left(  u\right)  >\sigma\left(  i\right)
\right]  +\left[  \sigma\left(  u\right)  >\sigma\left(  v\right)  \right] \\
&  \ \ \ \ \ \ \ \ \ \ +\sum_{i=u+1}^{v-1}\ell_{i}\left(  \tau\right)
+\sum_{i=u+1}^{v-1}\left(  \left[  \sigma\left(  i\right)  >\sigma\left(
v\right)  \right]  -\left[  \sigma\left(  i\right)  >\sigma\left(  u\right)
\right]  \right) \\
&  \ \ \ \ \ \ \ \ \ \ +\ell_{u}\left(  \tau\right)  -\sum_{i=u+1}%
^{v-1}\left[  \sigma\left(  v\right)  >\sigma\left(  i\right)  \right]
-\left[  \sigma\left(  v\right)  >\sigma\left(  u\right)  \right]
+\sum_{i=v+1}^{n}\ell_{i}\left(  \tau\right)
\end{align*}%
\begin{align*}
&  =\underbrace{\sum_{i=1}^{u-1}\ell_{i}\left(  \tau\right)  +\ell_{u}\left(
\tau\right)  +\sum_{i=u+1}^{v-1}\ell_{i}\left(  \tau\right)  +\ell_{v}\left(
\tau\right)  +\sum_{i=v+1}^{n}\ell_{i}\left(  \tau\right)  }_{\substack{=\ell
\left(  \tau\right)  \\\text{(by (\ref{pf.lem.sol.perm.lisitij.1.f.ltau=}))}%
}}\\
&  \ \ \ \ \ \ \ \ \ \ +\left[  \sigma\left(  u\right)  >\sigma\left(
v\right)  \right]  -\left[  \sigma\left(  v\right)  >\sigma\left(  u\right)
\right] \\
&  \ \ \ \ \ \ \ \ \ \ +\underbrace{\sum_{i=u+1}^{v-1}\left[  \sigma\left(
u\right)  >\sigma\left(  i\right)  \right]  +\sum_{i=u+1}^{v-1}\left(  \left[
\sigma\left(  i\right)  >\sigma\left(  v\right)  \right]  -\left[
\sigma\left(  i\right)  >\sigma\left(  u\right)  \right]  \right)
-\sum_{i=u+1}^{v-1}\left[  \sigma\left(  v\right)  >\sigma\left(  i\right)
\right]  }_{=\sum_{i=u+1}^{v-1}\left(  \left[  \sigma\left(  u\right)
>\sigma\left(  i\right)  \right]  +\left(  \left[  \sigma\left(  i\right)
>\sigma\left(  v\right)  \right]  -\left[  \sigma\left(  i\right)
>\sigma\left(  u\right)  \right]  \right)  -\left[  \sigma\left(  v\right)
>\sigma\left(  i\right)  \right]  \right)  }\\
&  =\ell\left(  \tau\right)  +\left[  \sigma\left(  u\right)  >\sigma\left(
v\right)  \right]  -\left[  \sigma\left(  v\right)  >\sigma\left(  u\right)
\right] \\
&  \ \ \ \ \ \ \ \ \ \ +\sum_{i=u+1}^{v-1}\underbrace{\left(  \left[
\sigma\left(  u\right)  >\sigma\left(  i\right)  \right]  +\left(  \left[
\sigma\left(  i\right)  >\sigma\left(  v\right)  \right]  -\left[
\sigma\left(  i\right)  >\sigma\left(  u\right)  \right]  \right)  -\left[
\sigma\left(  v\right)  >\sigma\left(  i\right)  \right]  \right)
}_{\substack{=\left[  \sigma\left(  u\right)  >\sigma\left(  i\right)
\right]  -\left[  \sigma\left(  v\right)  >\sigma\left(  i\right)  \right]
+\left[  \sigma\left(  i\right)  >\sigma\left(  v\right)  \right]  -\left[
\sigma\left(  i\right)  >\sigma\left(  u\right)  \right]  \\=2\left[
\sigma\left(  u\right)  >\sigma\left(  i\right)  >\sigma\left(  v\right)
\right]  -2\left[  \sigma\left(  v\right)  >\sigma\left(  i\right)
>\sigma\left(  u\right)  \right]  \\\text{(by Lemma
\ref{lem.sol.perm.lisitij.uiv})}}}\\
&  =\ell\left(  \tau\right)  +\left[  \sigma\left(  u\right)  >\sigma\left(
v\right)  \right]  -\left[  \sigma\left(  v\right)  >\sigma\left(  u\right)
\right] \\
&  \ \ \ \ \ \ \ \ \ \ +\sum_{i=u+1}^{v-1}\left(  2\left[  \sigma\left(
u\right)  >\sigma\left(  i\right)  >\sigma\left(  v\right)  \right]  -2\left[
\sigma\left(  v\right)  >\sigma\left(  i\right)  >\sigma\left(  u\right)
\right]  \right)  .
\end{align*}
This proves Lemma \ref{lem.sol.perm.lisitij.1} \textbf{(f)}.
\end{proof}

Our next lemma is precisely the claim of Exercise \ref{exe.perm.lisitij}, with
$i$, $j$ and $k$ renamed as $u$, $v$ and $i$:

\begin{lemma}
\label{lem.sol.perm.lisitij.2}Let $n\in\mathbb{N}$ and $\sigma\in S_{n}$. Let
$u$ and $v$ be two elements of $\left[  n\right]  $ such that $u<v$ and
$\sigma\left(  u\right)  >\sigma\left(  v\right)  $. Let $Q$ be the set of all
$i\in\left\{  u+1,u+2,\ldots,v-1\right\}  $ satisfying $\sigma\left(
u\right)  >\sigma\left(  i\right)  >\sigma\left(  v\right)  $. Then,%
\[
\ell\left(  \sigma\circ t_{u,v}\right)  =\ell\left(  \sigma\right)
-2\left\vert Q\right\vert -1.
\]

\end{lemma}

\begin{proof}
[Proof of Lemma \ref{lem.sol.perm.lisitij.2}.]Recall that $Q$ is the set of
all $i\in\left\{  u+1,u+2,\ldots,v-1\right\}  $ satisfying $\sigma\left(
u\right)  >\sigma\left(  i\right)  >\sigma\left(  v\right)  $. In other words,%
\begin{equation}
Q=\left\{  i\in\left\{  u+1,u+2,\ldots,v-1\right\}  \ \mid\ \sigma\left(
u\right)  >\sigma\left(  i\right)  >\sigma\left(  v\right)  \right\}  .
\label{pf.lem.sol.perm.lisitij.2.Q=}%
\end{equation}

Define a permutation $\tau\in S_{n}$ by $\tau=\sigma\circ t_{u,v}$. We have
$\sigma\left(  u\right)  >\sigma\left(  v\right)  $; thus, we don't have
$\sigma\left(  v\right)  >\sigma\left(  u\right)  $. Hence, $\left[
\sigma\left(  v\right)  >\sigma\left(  u\right)  \right]  =0$.

For any $i\in\left\{  u+1,u+2,\ldots,v-1\right\}  $, we have%
\begin{equation}
\left[  \sigma\left(  v\right)  >\sigma\left(  i\right)  >\sigma\left(
u\right)  \right]  =0 \label{pf.lem.sol.perm.lisitij.2.1}%
\end{equation}
(since we don't have $\sigma\left(  v\right)  >\sigma\left(  i\right)
>\sigma\left(  u\right)  $ (because we don't have $\sigma\left(  v\right)
>\sigma\left(  u\right)  $)).

Now, Lemma \ref{lem.sol.perm.lisitij.1} \textbf{(f)} yields%
\begin{align*}
\ell\left(  \sigma\right)   &  =\ell\left(  \tau\right)  +\underbrace{\left[
\sigma\left(  u\right)  >\sigma\left(  v\right)  \right]  }%
_{\substack{=1\\\text{(since }\sigma\left(  u\right)  >\sigma\left(  v\right)
\text{)}}}-\underbrace{\left[  \sigma\left(  v\right)  >\sigma\left(
u\right)  \right]  }_{=0}\\
&  \ \ \ \ \ \ \ \ \ \ +\underbrace{\sum_{i=u+1}^{v-1}}_{=\sum_{i\in\left\{
u+1,u+2,\ldots,v-1\right\}  }}\left(  2\left[  \sigma\left(  u\right)
>\sigma\left(  i\right)  >\sigma\left(  v\right)  \right]
-2\underbrace{\left[  \sigma\left(  v\right)  >\sigma\left(  i\right)
>\sigma\left(  u\right)  \right]  }_{\substack{=0\\\text{(by
(\ref{pf.lem.sol.perm.lisitij.2.1}))}}}\right) \\
&  =\ell\left(  \tau\right)  +\underbrace{1-0}_{=1}+\sum_{i\in\left\{
u+1,u+2,\ldots,v-1\right\}  }\underbrace{\left(  2\left[  \sigma\left(
u\right)  >\sigma\left(  i\right)  >\sigma\left(  v\right)  \right]
-2\cdot0\right)  }_{=2\left[  \sigma\left(  u\right)  >\sigma\left(  i\right)
>\sigma\left(  v\right)  \right]  }\\
&  =\ell\left(  \tau\right)  +1+\underbrace{\sum_{i\in\left\{  u+1,u+2,\ldots
,v-1\right\}  }2\left[  \sigma\left(  u\right)  >\sigma\left(  i\right)
>\sigma\left(  v\right)  \right]  }_{=2\sum_{i\in\left\{  u+1,u+2,\ldots
,v-1\right\}  }\left[  \sigma\left(  u\right)  >\sigma\left(  i\right)
>\sigma\left(  v\right)  \right]  }\\
&  =\ell\left(  \tau\right)  +1+2\sum_{i\in\left\{  u+1,u+2,\ldots
,v-1\right\}  }\left[  \sigma\left(  u\right)  >\sigma\left(  i\right)
>\sigma\left(  v\right)  \right]  .
\end{align*}
In view of%
\begin{align*}
&  \sum_{i\in\left\{  u+1,u+2,\ldots,v-1\right\}  }\left[  \sigma\left(
u\right)  >\sigma\left(  i\right)  >\sigma\left(  v\right)  \right] \\
&  =\sum_{\substack{i\in\left\{  u+1,u+2,\ldots,v-1\right\}  ;\\\sigma\left(
u\right)  >\sigma\left(  i\right)  >\sigma\left(  v\right)  }%
}\underbrace{\left[  \sigma\left(  u\right)  >\sigma\left(  i\right)
>\sigma\left(  v\right)  \right]  }_{\substack{=1\\\text{(since }\sigma\left(
u\right)  >\sigma\left(  i\right)  >\sigma\left(  v\right)  \text{)}}%
}+\sum_{\substack{i\in\left\{  u+1,u+2,\ldots,v-1\right\}  ;\\\text{not
}\sigma\left(  u\right)  >\sigma\left(  i\right)  >\sigma\left(  v\right)
}}\underbrace{\left[  \sigma\left(  u\right)  >\sigma\left(  i\right)
>\sigma\left(  v\right)  \right]  }_{\substack{=0\\\text{(since we don't have
}\sigma\left(  u\right)  >\sigma\left(  i\right)  >\sigma\left(  v\right)
\text{)}}}\\
&  \ \ \ \ \ \ \ \ \ \ \left(
\begin{array}
[c]{c}%
\text{since each }i\in\left\{  u+1,u+2,\ldots,v-1\right\}  \text{ satisfies
either }\left(  \sigma\left(  u\right)  >\sigma\left(  i\right)
>\sigma\left(  v\right)  \right) \\
\text{or }\left(  \text{not }\sigma\left(  u\right)  >\sigma\left(  i\right)
>\sigma\left(  v\right)  \right)  \text{ (but not both)}%
\end{array}
\right) \\
&  =\sum_{\substack{i\in\left\{  u+1,u+2,\ldots,v-1\right\}  ;\\\sigma\left(
u\right)  >\sigma\left(  i\right)  >\sigma\left(  v\right)  }%
}1+\underbrace{\sum_{\substack{i\in\left\{  u+1,u+2,\ldots,v-1\right\}
;\\\text{not }\sigma\left(  u\right)  >\sigma\left(  i\right)  >\sigma\left(
v\right)  }}0}_{=0}=\sum_{\substack{i\in\left\{  u+1,u+2,\ldots,v-1\right\}
;\\\sigma\left(  u\right)  >\sigma\left(  i\right)  >\sigma\left(  v\right)
}}1\\
&  =\left\vert \underbrace{\left\{  i\in\left\{  u+1,u+2,\ldots,v-1\right\}
\ \mid\ \sigma\left(  u\right)  >\sigma\left(  i\right)  >\sigma\left(
v\right)  \right\}  }_{\substack{=Q\\\text{(by
(\ref{pf.lem.sol.perm.lisitij.2.Q=}))}}}\right\vert \cdot1=\left\vert
Q\right\vert \cdot1=\left\vert Q\right\vert ,
\end{align*}
this becomes%
\[
\ell\left(  \sigma\right)  =\ell\left(  \tau\right)  +1+2\underbrace{\sum
_{i\in\left\{  u+1,u+2,\ldots,v-1\right\}  }\left[  \sigma\left(  u\right)
>\sigma\left(  i\right)  >\sigma\left(  v\right)  \right]  }_{=\left\vert
Q\right\vert }=\ell\left(  \tau\right)  +1+2\left\vert Q\right\vert .
\]
Solving this equation for $\ell\left(  \tau\right)  $, we find%
\[
\ell\left(  \tau\right)  =\ell\left(  \sigma\right)  -1-2\left\vert
Q\right\vert =\ell\left(  \sigma\right)  -2\left\vert Q\right\vert -1.
\]
In view of $\tau=\sigma\circ t_{u,v}$, this rewrites as $\ell\left(
\sigma\circ t_{u,v}\right)  =\ell\left(  \sigma\right)  -2\left\vert
Q\right\vert -1$. Thus, Lemma \ref{lem.sol.perm.lisitij.2} is proven.
\end{proof}

\begin{proof}
[Solution to Exercise \ref{exe.perm.lisitij}.]We obtain the statement of
Exercise \ref{exe.perm.lisitij} if we rename the variables $u$, $v$ and $i$ as
$i$, $j$ and $k$ in Lemma \ref{lem.sol.perm.lisitij.2}. Thus, Exercise
\ref{exe.perm.lisitij} is solved.
\end{proof}

\subsection{\label{sect.sol.perm.lehmer.rothe}Solution to Exercise
\ref{exe.perm.lehmer.rothe}}

Throughout Section \ref{sect.sol.perm.lehmer.rothe}, we shall use the same
notations that were used in Section \ref{sect.perm.lehmer}. We shall also use
Definition \ref{def.sol.perm.c=ttt.cuv}.

We begin with a simple lemma that will help us resolve parts \textbf{(a)} and
\textbf{(b)} of Exercise \ref{exe.perm.lehmer.rothe}:

\begin{lemma}
\label{lem.sol.perm.lehmer.ab}Let $n\in\mathbb{N}$ and $\sigma\in S_{n}$. Let
$i\in\left[  n\right]  $.

\textbf{(a)} The permutation $c_{i,i+\ell_{i}\left(  \sigma\right)  }\in
S_{n}$ is well-defined.

\textbf{(b)} We have $c_{i,i+\ell_{i}\left(  \sigma\right)  }=s_{i^{\prime}%
-1}\circ s_{i^{\prime}-2}\circ\cdots\circ s_{i}$, where $i^{\prime}=i+\ell
_{i}\left(  \sigma\right)  $.
\end{lemma}

\begin{proof}
[Proof of Lemma \ref{lem.sol.perm.lehmer.ab}.]We know (from Corollary
\ref{cor.sol.perm.c=ttt.cuv1} \textbf{(a)}) that the permutation $c_{u,v}$ is
well-defined whenever $u$ and $v$ are two elements of $\left[  n\right]  $
such that $u\leq v$.

Now, $i\in\left[  n\right]  =\left\{  1,2,\ldots,n\right\}  $ (by the
definition of $\left[  n\right]  $). Also, Proposition
\ref{prop.perm.lehmer.wd} yields%
\[
\left(  \ell_{1}\left(  \sigma\right)  ,\ell_{2}\left(  \sigma\right)
,\ldots,\ell_{n}\left(  \sigma\right)  \right)  \in H=\left[  n-1\right]
_{0}\times\left[  n-2\right]  _{0}\times\cdots\times\left[  n-n\right]  _{0}%
\]
(by the definition of $H$). In other words, $\ell_{j}\left(  \sigma\right)
\in\left[  n-j\right]  _{0}$ for each $j\in\left\{  1,2,\ldots,n\right\}  $.
Applying this to $j=i$, we obtain
\[
\ell_{i}\left(  \sigma\right)  \in\left[  n-i\right]  _{0}=\left\{
0,1,\ldots,n-i\right\}
\]
(by the definition of $\left[  n-i\right]  _{0}$). Hence, $\ell_{i}\left(
\sigma\right)  \geq0$ and $\ell_{i}\left(  \sigma\right)  \leq n-i$. Now,
$i+\underbrace{\ell_{i}\left(  \sigma\right)  }_{\geq0}\geq i$, so that $i\leq
i+\ell_{i}\left(  \sigma\right)  $.

\begin{vershort}
Also, $i+\underbrace{\ell_{i}\left(  \sigma\right)  }_{\geq0}\geq i\geq1$
(since $i\in\left\{  1,2,\ldots,n\right\}  $). Combining this with
$i+\underbrace{\ell_{i}\left(  \sigma\right)  }_{\leq n-i}\leq i+n-i=n$, we
obtain $i+\ell_{i}\left(  \sigma\right)  \in\left\{  1,2,\ldots,n\right\}
=\left[  n\right]  $. Hence, $i$ and $i+\ell_{i}\left(  \sigma\right)  $ are
two elements of $\left[  n\right]  $ such that $i\leq i+\ell_{i}\left(
\sigma\right)  $. Hence, Lemma \ref{cor.sol.perm.c=ttt.cuv1} \textbf{(a)}
(applied to $u=i$ and $v=i+\ell_{i}\left(  \sigma\right)  $) yields that the
permutation $c_{i,i+\ell_{i}\left(  \sigma\right)  }$ is well-defined. This
proves Lemma \ref{lem.sol.perm.lehmer.ab} \textbf{(a)}.
\end{vershort}

\begin{verlong}
Also, $i+\underbrace{\ell_{i}\left(  \sigma\right)  }_{\geq0}\geq i\geq1$
(since $i\in\left\{  1,2,\ldots,n\right\}  $). Combining this with
$i+\underbrace{\ell_{i}\left(  \sigma\right)  }_{\leq n-i}\leq i+n-i=n$, we
obtain $i+\ell_{i}\left(  \sigma\right)  \in\left\{  1,2,\ldots,n\right\}
=\left[  n\right]  $. Hence, $i$ and $i+\ell_{i}\left(  \sigma\right)  $ are
two elements of $\left[  n\right]  $ (since $i\in\left[  n\right]  $ and
$i+\ell_{i}\left(  \sigma\right)  \in\left[  n\right]  $) such that $i\leq
i+\ell_{i}\left(  \sigma\right)  $. Hence, Lemma \ref{cor.sol.perm.c=ttt.cuv1}
\textbf{(a)} (applied to $u=i$ and $v=i+\ell_{i}\left(  \sigma\right)  $)
yields that the permutation $c_{i,i+\ell_{i}\left(  \sigma\right)  }$ is
well-defined. This proves Lemma \ref{lem.sol.perm.lehmer.ab} \textbf{(a)}.
\end{verlong}

\textbf{(b)} Let $i^{\prime}=i+\ell_{i}\left(  \sigma\right)  $. We must prove
that $c_{i,i+\ell_{i}\left(  \sigma\right)  }=s_{i^{\prime}-1}\circ
s_{i^{\prime}-2}\circ\cdots\circ s_{i}$.

We have $i^{\prime}=i+\ell_{i}\left(  \sigma\right)  \geq i$, so that $i\leq
i^{\prime}$. Also, $i\in\left[  n\right]  $ and $i^{\prime}=i+\ell_{i}\left(
\sigma\right)  \in\left[  n\right]  $. Hence, Lemma
\ref{cor.sol.perm.c=ttt.cuv1} \textbf{(b)} (applied to $u=i$ and $v=i^{\prime
}$) yields $c_{i,i^{\prime}}=s_{i^{\prime}-1}\circ s_{i^{\prime}-2}\circ
\cdots\circ s_{i}$. But $i+\ell_{i}\left(  \sigma\right)  =i^{\prime}$ (since
$i^{\prime}=i+\ell_{i}\left(  \sigma\right)  $); thus, $c_{i,i+\ell_{i}\left(
\sigma\right)  }=c_{i,i^{\prime}}=s_{i^{\prime}-1}\circ s_{i^{\prime}-2}%
\circ\cdots\circ s_{i}$. This proves Lemma \ref{lem.sol.perm.lehmer.ab}
\textbf{(b)}.
\end{proof}

Next, we state a property of the permutations $c_{u,v}$ that is essentially
clear from their definition:

\begin{lemma}
\label{lem.sol.perm.lehmer.cuv.values}Let $n\in\mathbb{N}$. Let $u$ and $v$ be
two elements of $\left[  n\right]  $ such that $u\leq v$. Consider the
permutation $c_{u,v}\in S_{n}$ defined in Definition
\ref{def.sol.perm.c=ttt.cuv}. Then:

\textbf{(a)} We have $c_{u,v}\left(  q\right)  =q$ for each $q\in\left[
n\right]  \setminus\left\{  u,u+1,\ldots,v\right\}  $.

\textbf{(b)} We have $c_{u,v}\left(  q\right)  =q-1$ for each $q\in\left\{
u+1,u+2,\ldots,v\right\}  $.

\textbf{(c)} We have $c_{u,v}\left(  u\right)  =v$.
\end{lemma}

\begin{vershort}
\begin{proof}
[Proof of Lemma \ref{lem.sol.perm.lehmer.cuv.values}.]The definition of
$c_{u,v}$ yields $c_{u,v}=\operatorname*{cyc}\nolimits_{v,v-1,v-2,\ldots,u}$.
But the $\left(  v-u+1\right)  $-cycle $\operatorname*{cyc}%
\nolimits_{v,v-1,v-2,\ldots,u}$ was defined as the permutation in $S_{n}$
which sends $v,v-1,\ldots,u+1,u$ to $v-1,v-2,\ldots,u,v$, respectively, while
leaving all other elements of $\left[  n\right]  $ fixed. In other words,
$c_{u,v}$ is the permutation in $S_{n}$ which sends $v,v-1,\ldots,u+1,u$ to
$v-1,v-2,\ldots,u,v$, respectively, while leaving all other elements of
$\left[  n\right]  $ fixed. This immediately yields all three parts of Lemma
\ref{lem.sol.perm.lehmer.cuv.values}.
\end{proof}
\end{vershort}

\begin{verlong}
\begin{proof}
[Proof of Lemma \ref{lem.sol.perm.lehmer.cuv.values}.]The definition of
$c_{u,v}$ yields $c_{u,v}=\operatorname*{cyc}\nolimits_{v,v-1,v-2,\ldots,u}$.

We have $u\in\left[  n\right]  =\left\{  1,2,\ldots,n\right\}  $ (by the
definition of $\left[  n\right]  $), thus $u\geq1$. Hence, $1\leq u$. Also,
$v\in\left[  n\right]  =\left\{  1,2,\ldots,n\right\}  $, thus $v\leq n$.
Hence, $1\leq u\leq v\leq n$.

Now, $\underbrace{v}_{\leq n}-\underbrace{u}_{\geq1}\leq n-1$. Combining this
with $v-\underbrace{u}_{\leq v}\geq v-v=0$, we obtain $v-u\in\left\{
0,1,\ldots,n-1\right\}  $, so that $\left(  v-u\right)  +1\in\left\{
1,2,\ldots,n\right\}  $. Hence, $v-u+1=\left(  v-u\right)  +1\in\left\{
1,2,\ldots,n\right\}  $.

Let $k=v-u+1$. Then, $k=v-u+1\in\left\{  1,2,\ldots,n\right\}  $. Also,
$v-\underbrace{k}_{=v-u+1}+1=v-\left(  v-u+1\right)  +1=u$ and $v-u=k-1$
(since $v-u+1=k$).

For each $p\in\left\{  1,2,\ldots,k\right\}  $, define the integer $i_{p}$ by
$i_{p}=v-p+1$. Thus,%
\begin{align*}
\left(  i_{1},i_{2},\ldots,i_{k}\right)   &  =\left(  v-1+1,v-2+1,\ldots
,v-k+1\right) \\
&  =\left(  v,v-1,v-2,\ldots,v-k+1\right)  =\left(  v,v-1,v-2,\ldots,u\right)
\end{align*}
(since $v-k+1=u$). In other words, $\left(  v,v-1,v-2,\ldots,u\right)
=\left(  i_{1},i_{2},\ldots,i_{k}\right)  $. Thus, $\operatorname*{cyc}%
\nolimits_{v,v-1,v-2,\ldots,u}=\operatorname*{cyc}\nolimits_{i_{1}%
,i_{2},\ldots,i_{k}}$. Hence,%
\begin{equation}
c_{u,v}=\operatorname*{cyc}\nolimits_{v,v-1,v-2,\ldots,u}=\operatorname*{cyc}%
\nolimits_{i_{1},i_{2},\ldots,i_{k}}.
\label{pf.lem.sol.perm.lehmer.cuv.values.cuv=}%
\end{equation}

We have defined $\operatorname*{cyc}\nolimits_{i_{1},i_{2},\ldots,i_{k}}$ to
be the permutation in $S_{n}$ which sends $i_{1},i_{2},\ldots,i_{k}$ to
$i_{2},i_{3},\ldots,i_{k},i_{1}$, respectively, while leaving all other
elements of $\left[  n\right]  $ fixed. Therefore:

\begin{itemize}
\item The permutation $\operatorname*{cyc}\nolimits_{i_{1},i_{2},\ldots,i_{k}%
}$ sends $i_{1},i_{2},\ldots,i_{k}$ to $i_{2},i_{3},\ldots,i_{k},i_{1}$,
respectively. In other words,%
\begin{equation}
\operatorname*{cyc}\nolimits_{i_{1},i_{2},\ldots,i_{k}}\left(  i_{p}\right)
=i_{p+1}\ \ \ \ \ \ \ \ \ \ \text{for every }p\in\left\{  1,2,\ldots
,k\right\}  , \label{pf.lem.sol.perm.lehmer.cuv.values.cyc-manifest.1}%
\end{equation}
where $i_{k+1}$ means $i_{1}$.

\item The permutation $\operatorname*{cyc}\nolimits_{i_{1},i_{2},\ldots,i_{k}%
}$ leaves all other elements of $\left[  n\right]  $ fixed (where
\textquotedblleft other\textquotedblright\ means \textquotedblleft other than
$i_{1},i_{2},\ldots,i_{k}$\textquotedblright). In other words,%
\begin{equation}
\operatorname*{cyc}\nolimits_{i_{1},i_{2},\ldots,i_{k}}\left(  q\right)
=q\ \ \ \ \ \ \ \ \ \ \text{for every }q\in\left[  n\right]  \setminus\left\{
i_{1},i_{2},\ldots,i_{k}\right\}  .
\label{pf.lem.sol.perm.lehmer.cuv.values.cyc-manifest.2}%
\end{equation}

\end{itemize}

Let us define $i_{k+1}$ to be a synonym for $i_{1}$. Thus, $i_{k+1}=i_{1}$.

\textbf{(a)} Let $q\in\left[  n\right]  \setminus\left\{  u,u+1,\ldots
,v\right\}  $. Now,%
\[
\left\{  u,u+1,\ldots,v\right\}  =\left\{  v,v-1,v-2,\ldots,u\right\}
=\left\{  i_{1},i_{2},\ldots,i_{k}\right\}
\]
(since $\left(  v,v-1,v-2,\ldots,u\right)  =\left(  i_{1},i_{2},\ldots
,i_{k}\right)  $). Hence, $q\in\left[  n\right]  \setminus\underbrace{\left\{
u,u+1,\ldots,v\right\}  }_{=\left\{  i_{1},i_{2},\ldots,i_{k}\right\}
}=\left[  n\right]  \setminus\left\{  i_{1},i_{2},\ldots,i_{k}\right\}  $.
Thus, (\ref{pf.lem.sol.perm.lehmer.cuv.values.cyc-manifest.2}) yields
$\operatorname*{cyc}\nolimits_{i_{1},i_{2},\ldots,i_{k}}\left(  q\right)  =q$.
Now, $\underbrace{c_{u,v}}_{\substack{=\operatorname*{cyc}\nolimits_{i_{1}%
,i_{2},\ldots,i_{k}}\\\text{(by (\ref{pf.lem.sol.perm.lehmer.cuv.values.cuv=}%
))}}}\left(  q\right)  =\operatorname*{cyc}\nolimits_{i_{1},i_{2},\ldots
,i_{k}}\left(  q\right)  =q$. This proves Lemma
\ref{lem.sol.perm.lehmer.cuv.values} \textbf{(a)}.

\textbf{(b)} Let $q\in\left\{  u+1,u+2,\ldots,v\right\}  $. Thus,
\begin{align*}
q  &  \in\left\{  u+1,u+2,\ldots,v\right\}  \subseteq\left\{  1+1,1+2,\ldots
,n\right\}  \ \ \ \ \ \ \ \ \ \ \left(  \text{since }u\geq1\text{ and }v\leq
n\right) \\
&  =\left\{  2,3,\ldots,n\right\}  \subseteq\left\{  1,2,\ldots,n\right\}
=\left[  n\right]  .
\end{align*}
Also, define an integer $p$ by $p=v-q+1$. Thus,
\begin{align*}
p  &  =v-q+1=\left(  v+1\right)  -q\in\left\{  1,2,\ldots,v-u\right\}
\ \ \ \ \ \ \ \ \ \ \left(  \text{since }q\in\left\{  u+1,u+2,\ldots
,v\right\}  \right) \\
&  =\left\{  1,2,\ldots,k-1\right\}  \ \ \ \ \ \ \ \ \ \ \left(  \text{since
}v-u=k-1\right) \\
&  \subseteq\left\{  1,2,\ldots,k\right\}  .
\end{align*}
Hence, (\ref{pf.lem.sol.perm.lehmer.cuv.values.cyc-manifest.1}) yields
$\operatorname*{cyc}\nolimits_{i_{1},i_{2},\ldots,i_{k}}\left(  i_{p}\right)
=i_{p+1}$. Now, $\underbrace{c_{u,v}}_{\substack{=\operatorname*{cyc}%
\nolimits_{i_{1},i_{2},\ldots,i_{k}}\\\text{(by
(\ref{pf.lem.sol.perm.lehmer.cuv.values.cuv=}))}}}\left(  i_{p}\right)
=\operatorname*{cyc}\nolimits_{i_{1},i_{2},\ldots,i_{k}}\left(  i_{p}\right)
=i_{p+1}$.

But $p\in\left\{  1,2,\ldots,k\right\}  $; thus, the definition of $i_{p}$
yields $i_{p}=v-\underbrace{p}_{=v-q+1}+1=v-\left(  v-q+1\right)  +1=q$. Also,
$p\in\left\{  1,2,\ldots,k-1\right\}  $, so that $p+1\in\left\{
2,3,\ldots,k\right\}  \subseteq\left\{  1,2,\ldots,k\right\}  $; hence, the
definition of $i_{p+1}$ yields $i_{p+1}=v-\left(  p+1\right)
+1=v-\underbrace{p}_{=v-q+1}=v-\left(  v-q+1\right)  =q-1$.

But we have shown that $c_{u,v}\left(  i_{p}\right)  =i_{p+1}$. In view of
$i_{p}=q$ and $i_{p+1}=q-1$, this rewrites as $c_{u,v}\left(  q\right)  =q-1$.
This proves Lemma \ref{lem.sol.perm.lehmer.cuv.values} \textbf{(b)}.

\textbf{(c)} We have $u\in\left[  n\right]  $. Also, $k\geq1$ (since
$k\in\left\{  1,2,\ldots,n\right\}  $) and thus $k\in\left\{  1,2,\ldots
,k\right\}  $. Hence, (\ref{pf.lem.sol.perm.lehmer.cuv.values.cyc-manifest.1})
(applied to $p=k$) yields $\operatorname*{cyc}\nolimits_{i_{1},i_{2}%
,\ldots,i_{k}}\left(  i_{k}\right)  =i_{k+1}=i_{1}$. Now, $\underbrace{c_{u,v}%
}_{\substack{=\operatorname*{cyc}\nolimits_{i_{1},i_{2},\ldots,i_{k}%
}\\\text{(by (\ref{pf.lem.sol.perm.lehmer.cuv.values.cuv=}))}}}\left(
i_{k}\right)  =\operatorname*{cyc}\nolimits_{i_{1},i_{2},\ldots,i_{k}}\left(
i_{k}\right)  =i_{1}$.

But $k\in\left\{  1,2,\ldots,k\right\}  $; thus, the definition of $i_{k}$
yields $i_{k}=v-k+1=u$. Also, $1\in\left\{  1,2,\ldots,k\right\}  $ (since
$k\geq1$); hence, the definition of $i_{1}$ yields $i_{1}=v-1+1=v$.

But we have shown that $c_{u,v}\left(  i_{k}\right)  =i_{1}$. In view of
$i_{k}=u$ and $i_{1}=v$, this rewrites as $c_{u,v}\left(  u\right)  =v$. This
proves Lemma \ref{lem.sol.perm.lehmer.cuv.values} \textbf{(c)}.
\end{proof}
\end{verlong}

\begin{corollary}
\label{cor.sol.perm.lehmer.cuv.notfar}Let $n\in\mathbb{N}$. Let $u$ and $v$ be
two elements of $\left[  n\right]  $ such that $u\leq v$. Consider the
permutation $c_{u,v}\in S_{n}$ defined in Definition
\ref{def.sol.perm.c=ttt.cuv}. Let $k\in\left[  n\right]  $ be such that $k\neq
u$. Then, $c_{u,v}\left(  k\right)  \in\left\{  k-1,k\right\}  $.
\end{corollary}

\begin{proof}
[Proof of Corollary \ref{cor.sol.perm.lehmer.cuv.notfar}.]We are in one of the
following two cases:

\textit{Case 1:} We have $k\in\left\{  u,u+1,\ldots,v\right\}  $.

\textit{Case 2:} We don't have $k\in\left\{  u,u+1,\ldots,v\right\}  $.

Let us first consider Case 1. In this case, we have $k\in\left\{
u,u+1,\ldots,v\right\}  $. Combining this with $k\neq u$, we obtain
$k\in\left\{  u,u+1,\ldots,v\right\}  \setminus\left\{  u\right\}  =\left\{
u+1,u+2,\ldots,v\right\}  $. Hence, Lemma \ref{lem.sol.perm.lehmer.cuv.values}
\textbf{(b)} (applied to $q=k$) yields $c_{u,v}\left(  k\right)
=k-1\in\left\{  k-1,k\right\}  $. Thus, Corollary
\ref{cor.sol.perm.lehmer.cuv.notfar} is proven in Case 1.

Let us now consider Case 2. In this case, we don't have $k\in\left\{
u,u+1,\ldots,v\right\}  $. In other words, we have $k\notin\left\{
u,u+1,\ldots,v\right\}  $. Combining $k\in\left[  n\right]  $ with
$k\notin\left\{  u,u+1,\ldots,v\right\}  $, we obtain $k\in\left[  n\right]
\setminus\left\{  u,u+1,\ldots,v\right\}  $. Hence, Lemma
\ref{lem.sol.perm.lehmer.cuv.values} \textbf{(a)} (applied to $q=k$) yields
$c_{u,v}\left(  k\right)  =k\in\left\{  k-1,k\right\}  $. Thus, Corollary
\ref{cor.sol.perm.lehmer.cuv.notfar} is proven in Case 2.

We have now proven Corollary \ref{cor.sol.perm.lehmer.cuv.notfar} in each of
the two Cases 1 and 2. Since these two Cases cover all possibilities, we thus
conclude that Corollary \ref{cor.sol.perm.lehmer.cuv.notfar} always holds.
\end{proof}

Next comes another simple property of the Lehmer code:

\begin{lemma}
\label{lem.sol.perm.lehmer.lu-if-fix}Let $n\in\mathbb{N}$. Let $u\in\left[
n\right]  $. Let $\sigma\in S_{n}$. Assume that
\begin{equation}
\left(  \sigma\left(  q\right)  =q\text{ for each }q\in\left\{  1,2,\ldots
,u-1\right\}  \right)  . \label{eq.lem.sol.perm.lehmer.lu-if-fix.ass}%
\end{equation}
Then, $u+\ell_{u}\left(  \sigma\right)  =\sigma\left(  u\right)  $.
\end{lemma}

\begin{proof}
[Proof of Lemma \ref{lem.sol.perm.lehmer.lu-if-fix}.]We have%
\begin{align}
\sigma\left(  \left\{  1,2,\ldots,u-1\right\}  \right)   &  =\left\{
\underbrace{\sigma\left(  q\right)  }_{\substack{=q\\\text{(by
(\ref{eq.lem.sol.perm.lehmer.lu-if-fix.ass}))}}}\ \mid\ q\in\left\{
1,2,\ldots,u-1\right\}  \right\} \nonumber\\
&  =\left\{  q\ \mid\ q\in\left\{  1,2,\ldots,u-1\right\}  \right\}
\nonumber\\
&  =\left\{  1,2,\ldots,u-1\right\}  .
\label{pf.lem.sol.perm.lehmer.lu-if-fix.1}%
\end{align}

\begin{verlong}
We have $\left[  n\right]  =\left\{  1,2,\ldots,n\right\}  $ (by the
definition of $\left[  n\right]  $).

Recall that $S_{n}$ is the set of all permutations of the set $\left\{
1,2,\ldots,n\right\}  $. In other words, $S_{n}$ is the set of all
permutations of the set $\left[  n\right]  $ (since $\left[  n\right]
=\left\{  1,2,\ldots,n\right\}  $). We have $\sigma\in S_{n}$. In other words,
$\sigma$ is a permutation of $\left[  n\right]  $ (since $S_{n}$ is the set of
all permutations of the set $\left[  n\right]  $). In other words, $\sigma$ is
a bijective map $\left[  n\right]  \rightarrow\left[  n\right]  $. Hence, the
map $\sigma$ is bijective, and thus injective.

We know that $\sigma$ is a map $\left[  n\right]  \rightarrow\left[  n\right]
$. Thus, $\sigma\left(  u\right)  \in\left[  n\right]  =\left\{
1,2,\ldots,n\right\}  $, so that $\sigma\left(  u\right)  \geq1$.
\end{verlong}

Next, we claim that $\sigma\left(  u\right)  \geq u$.

\begin{vershort}
[\textit{Proof:} Assume the contrary. Thus, $\sigma\left(  u\right)  <u$.
Hence, $\sigma\left(  u\right)  \in\left\{  1,2,\ldots,u-1\right\}  $ (because
$\sigma\left(  u\right)  $ is a positive integer). Thus,%
\[
\sigma\left(  u\right)  \in\left\{  1,2,\ldots,u-1\right\}  =\sigma\left(
\left\{  1,2,\ldots,u-1\right\}  \right)  \ \ \ \ \ \ \ \ \ \ \left(  \text{by
(\ref{pf.lem.sol.perm.lehmer.lu-if-fix.1})}\right)  .
\]
In other words, there exists some $q\in\left\{  1,2,\ldots,u-1\right\}  $ such
that $\sigma\left(  u\right)  =\sigma\left(  q\right)  $. Consider this $q$.
But $\sigma\in S_{n}$; thus, $\sigma$ is a permutation of $\left\{
1,2,\ldots,n\right\}  $ (since $S_{n}$ is the set of all permutations of
$\left\{  1,2,\ldots,n\right\}  $). Hence, $\sigma$ is a bijective map, thus
an injective map. Thus, from $\sigma\left(  u\right)  =\sigma\left(  q\right)
$, we obtain $u=q\in\left\{  1,2,\ldots,u-1\right\}  $, whence $u\leq u-1<u$.
This is absurd. Hence, we have obtained a contradiction. Thus, our assumption
was false. Thus, $\sigma\left(  u\right)  \geq u$ is proven.]
\end{vershort}

\begin{verlong}
[\textit{Proof:} Assume the contrary. Thus, $\sigma\left(  u\right)  <u$.
Hence, $\sigma\left(  u\right)  \leq u-1$ (since $\sigma\left(  u\right)  $
and $u$ are integers).

Combining this with $\sigma\left(  u\right)  \geq1$, we obtain $\sigma\left(
u\right)  \in\left\{  1,2,\ldots,u-1\right\}  $. Hence,
(\ref{eq.lem.sol.perm.lehmer.lu-if-fix.ass}) (applied to $q=\sigma\left(
u\right)  $) yields $\sigma\left(  \sigma\left(  u\right)  \right)
=\sigma\left(  u\right)  $. Since $\sigma$ is injective, we thus conclude
$\sigma\left(  u\right)  =u$. This contradicts $\sigma\left(  u\right)  <u$.

Hence, we have obtained a contradiction. Thus, our assumption was false. Thus,
$\sigma\left(  u\right)  \geq u$ is proven.]
\end{verlong}

If $a$ and $b$ are any two nonnegative integers satisfying $a\geq b$, then%
\[
\left\{  1,2,\ldots,a\right\}  \setminus\left\{  1,2,\ldots,b\right\}
=\left\{  b+1,b+2,\ldots,a\right\}
\]
and thus%
\begin{align}
\left\vert \left\{  1,2,\ldots,a\right\}  \setminus\left\{  1,2,\ldots
,b\right\}  \right\vert  &  =\left\vert \left\{  b+1,b+2,\ldots,a\right\}
\right\vert \nonumber\\
&  =a-b\ \ \ \ \ \ \ \ \ \ \left(  \text{since }a\geq b\right)  .
\label{pf.lem.sol.perm.lehmer.lu-if-fix.ab}%
\end{align}

\begin{vershort}
But $\sigma\left(  u\right)  \geq1$ (since $\sigma\left(  u\right)  \in\left[
n\right]  $) and $u\geq1$ (since $u\in\left[  n\right]  $). Hence,
$\sigma\left(  u\right)  -1$ and $u-1$ are nonnegative integers. These
integers satisfy $\sigma\left(  u\right)  -1\geq u-1$ (since $\sigma\left(
u\right)  \geq u$).
\end{vershort}

\begin{verlong}
We have proven that $\sigma\left(  u\right)  \geq u$. Thus, $\sigma\left(
u\right)  -1\geq u-1$. Also, $u\in\left[  n\right]  =\left\{  1,2,\ldots
,n\right\}  $ (by the definition of $\left[  n\right]  $), so that $u\geq1$
and thus $u-1\geq0$; hence, $u-1$ is a nonnegative integer. Also,
$\sigma\left(  u\right)  \geq1$, thus $\sigma\left(  u\right)  -1\geq0$, so
that $\sigma\left(  u\right)  -1$ is a nonnegative integer. Thus,
$\sigma\left(  u\right)  -1$ and $u-1$ are nonnegative integers satisfying
$\sigma\left(  u\right)  -1\geq u-1$.
\end{verlong}

Hence, (\ref{pf.lem.sol.perm.lehmer.lu-if-fix.ab}) (applied to $a=\sigma
\left(  u\right)  -1$ and $b=u-1$) yields%
\[
\left\vert \left\{  1,2,\ldots,\sigma\left(  u\right)  -1\right\}
\setminus\left\{  1,2,\ldots,u-1\right\}  \right\vert =\left(  \sigma\left(
u\right)  -1\right)  -\left(  u-1\right)  =\sigma\left(  u\right)  -u.
\]
But Lemma \ref{lem.perm.lexico1.lis} \textbf{(b)} (applied to $i=u$) yields
\begin{align*}
\ell_{u}\left(  \sigma\right)   &  =\left\vert \underbrace{\left[
\sigma\left(  u\right)  -1\right]  }_{\substack{=\left\{  1,2,\ldots
,\sigma\left(  u\right)  -1\right\}  \\\text{(by the definition of }\left[
\sigma\left(  u\right)  -1\right]  \text{)}}}\setminus\sigma\left(
\underbrace{\left[  u-1\right]  }_{\substack{=\left\{  1,2,\ldots,u-1\right\}
\\\text{(by the definition of }\left[  u-1\right]  \text{)}}}\right)
\right\vert \\
&  =\left\vert \left\{  1,2,\ldots,\sigma\left(  u\right)  -1\right\}
\setminus\underbrace{\sigma\left(  \left\{  1,2,\ldots,u-1\right\}  \right)
}_{\substack{=\left\{  1,2,\ldots,u-1\right\}  \\\text{(by
(\ref{pf.lem.sol.perm.lehmer.lu-if-fix.1}))}}}\right\vert \\
&  =\left\vert \left\{  1,2,\ldots,\sigma\left(  u\right)  -1\right\}
\setminus\left\{  1,2,\ldots,u-1\right\}  \right\vert \\
&  =\sigma\left(  u\right)  -u.
\end{align*}
In other words, $u+\ell_{u}\left(  \sigma\right)  =\sigma\left(  u\right)  $.
This proves Lemma \ref{lem.sol.perm.lehmer.lu-if-fix}.
\end{proof}

Next, we state a lemma that will be crucial for our solution:

\begin{lemma}
\label{lem.sol.perm.lehmer.lili}Let $n\in\mathbb{N}$. Let $u$ and $v$ be two
elements of $\left[  n\right]  $ such that $u\leq v$. Let $\sigma\in S_{n}$ be
such that $\sigma\left(  u\right)  =v$. Let $\tau\in S_{n}$ be such that
$\tau=\left(  c_{u,v}\right)  ^{-1}\circ\sigma$. Then:

\textbf{(a)} We have $\tau\left(  u\right)  =u$.

\textbf{(b)} We have $\ell_{i}\left(  \tau\right)  =\ell_{i}\left(
\sigma\right)  $ for each $i\in\left\{  u+1,u+2,\ldots,n\right\}  $.

\textbf{(c)} If we have%
\begin{equation}
\left(  \sigma\left(  q\right)  =q\text{ for each }q\in\left\{  1,2,\ldots
,u-1\right\}  \right)  , \label{eq.lem.sol.perm.lehmer.lili.c.ass}%
\end{equation}
then we have%
\begin{equation}
\left(  \tau\left(  q\right)  =q\text{ for each }q\in\left\{  1,2,\ldots
,u\right\}  \right)  . \label{eq.lem.sol.perm.lehmer.lili.c.claim}%
\end{equation}

\end{lemma}

\begin{proof}
[Proof of Lemma \ref{lem.sol.perm.lehmer.lili}.]We have%
\begin{equation}
c_{u,v}\circ\underbrace{\tau}_{=\left(  c_{u,v}\right)  ^{-1}\circ\sigma
}=\underbrace{c_{u,v}\circ\left(  c_{u,v}\right)  ^{-1}}_{=\operatorname*{id}%
}\circ\sigma=\sigma. \label{pf.lem.sol.perm.lehmer.lili.rev}%
\end{equation}

Also, Lemma \ref{lem.sol.perm.lehmer.cuv.values} \textbf{(c)} yields
$c_{u,v}\left(  u\right)  =v$. Hence, $\left(  c_{u,v}\right)  ^{-1}\left(
v\right)  =u$. Now,%
\[
\underbrace{\tau}_{=\left(  c_{u,v}\right)  ^{-1}\circ\sigma}\left(  u\right)
=\left(  \left(  c_{u,v}\right)  ^{-1}\circ\sigma\right)  \left(  u\right)
=\left(  c_{u,v}\right)  ^{-1}\left(  \underbrace{\sigma\left(  u\right)
}_{=v}\right)  =\left(  c_{u,v}\right)  ^{-1}\left(  v\right)  =u.
\]
This proves Lemma \ref{lem.sol.perm.lehmer.lili} \textbf{(a)}.

\begin{verlong}
We have $\left[  n\right]  =\left\{  1,2,\ldots,n\right\}  $ (by the
definition of $\left[  n\right]  $).
\end{verlong}

Recall that $S_{n}$ is the set of all permutations of $\left\{  1,2,\ldots
,n\right\}  $. In other words, $S_{n}$ is the set of all permutations of
$\left[  n\right]  $ (since $\left[  n\right]  =\left\{  1,2,\ldots,n\right\}
$). We have $\sigma\in S_{n}$. In other words, $\sigma$ is a permutation of
$\left[  n\right]  $ (since $S_{n}$ is the set of all permutations of $\left[
n\right]  $). In other words, $\sigma$ is a bijective map $\left[  n\right]
\rightarrow\left[  n\right]  $. The same argument (applied to $\tau$ instead
of $\sigma$) shows that $\tau$ is a bijective map $\left[  n\right]
\rightarrow\left[  n\right]  $. The maps $\sigma$ and $\tau$ are bijective,
and thus injective.

\textbf{(b)} First, let us prove the following auxiliary claim:

\begin{statement}
\textit{Claim 1:} Let $k\in\left[  n\right]  $ be such that $k\neq u$. Then,
$\tau\left(  k\right)  -1\leq\sigma\left(  k\right)  \leq\tau\left(  k\right)
$.
\end{statement}

[\textit{Proof of Claim 1:} If we had $\tau\left(  k\right)  =\tau\left(
u\right)  $, then we would have $k=u$ (since $\tau$ is injective), which would
contradict $k\neq u$. Thus, we cannot have $\tau\left(  k\right)  =\tau\left(
u\right)  $. Hence, we have $\tau\left(  k\right)  \neq\tau\left(  u\right)
=u$ (by Lemma \ref{lem.sol.perm.lehmer.lili} \textbf{(a)}).

\begin{verlong}
We have $\tau\left(  k\right)  \in\left[  n\right]  $ (since $\tau$ is a map
$\left[  n\right]  \rightarrow\left[  n\right]  $) and $\tau\left(  k\right)
\neq u$.
\end{verlong}

Hence, Corollary \ref{cor.sol.perm.lehmer.cuv.notfar} (applied to $\tau\left(
k\right)  $ instead of $k$) yields \newline$c_{u,v}\left(  \tau\left(
k\right)  \right)  \in\left\{  \tau\left(  k\right)  -1,\tau\left(  k\right)
\right\}  $. Now,
\begin{align*}
\underbrace{\sigma}_{\substack{=c_{u,v}\circ\tau\\\text{(by
(\ref{pf.lem.sol.perm.lehmer.lili.rev}))}}}\left(  k\right)   &  =\left(
c_{u,v}\circ\tau\right)  \left(  k\right)  =c_{u,v}\left(  \tau\left(
k\right)  \right)  \in\left\{  \tau\left(  k\right)  -1,\tau\left(  k\right)
\right\} \\
&  =\left\{  \tau\left(  k\right)  -1,\tau\left(  k\right)  ,\ldots
,\tau\left(  k\right)  \right\}  .
\end{align*}
Hence, $\tau\left(  k\right)  -1\leq\sigma\left(  k\right)  \leq\tau\left(
k\right)  $. This proves Claim 1.]

Now, let $i\in\left\{  u+1,u+2,\ldots,n\right\}  $. We shall now show the
following claim:

\begin{statement}
\textit{Claim 2:} Let $j\in\left\{  i+1,i+2,\ldots,n\right\}  $. Then, we have
the following logical equivalence:%
\[
\left(  \tau\left(  i\right)  >\tau\left(  j\right)  \right)
\ \Longleftrightarrow\ \left(  \sigma\left(  i\right)  >\sigma\left(
j\right)  \right)  .
\]

\end{statement}

[\textit{Proof of Claim 2:} We have $j\in\left\{  i+1,i+2,\ldots,n\right\}  $,
thus $j\geq i+1>i$, and therefore $j\neq i$.

We have $i\in\left\{  u+1,u+2,\ldots,n\right\}  $, thus $i\geq u+1>u$, and
thus $i\neq u$. Now, $j>i>u$, thus $j\neq u$.

\begin{vershort}
Also, $i\in\left\{  u+1,u+2,\ldots,n\right\}  \subseteq\left[  n\right]  $ and
$j\in\left\{  i+1,i+2,\ldots,n\right\}  \subseteq\left[  n\right]  $ (since
$i\in\left[  n\right]  $).
\end{vershort}

\begin{verlong}
Also, $u\in\left[  n\right]  =\left\{  1,2,\ldots,n\right\}  $, so that
$u\geq1$ and therefore $u+1\geq u\geq1$. Now,
\begin{align*}
i  &  \in\left\{  u+1,u+2,\ldots,n\right\}  \subseteq\left\{  1,2,\ldots
,n\right\}  \ \ \ \ \ \ \ \ \ \ \left(  \text{since }u+1\geq1\right) \\
&  =\left[  n\right]  .
\end{align*}
Also, $i\in\left[  n\right]  =\left\{  1,2,\ldots,n\right\}  $, so that
$i\geq1$ and therefore $i+1\geq i\geq1$. Now,
\begin{align*}
j  &  \in\left\{  i+1,i+2,\ldots,n\right\}  \subseteq\left\{  1,2,\ldots
,n\right\}  \ \ \ \ \ \ \ \ \ \ \left(  \text{since }i+1\geq1\right) \\
&  =\left[  n\right]  .
\end{align*}

\end{verlong}

Claim 1 (applied to $k=i$) yields $\tau\left(  i\right)  -1\leq\sigma\left(
i\right)  \leq\tau\left(  i\right)  $ (since $i\in\left[  n\right]  $ and
$i\neq u$).

Claim 1 (applied to $k=j$) yields $\tau\left(  j\right)  -1\leq\sigma\left(
j\right)  \leq\tau\left(  j\right)  $ (since $j\in\left[  n\right]  $ and
$j\neq u$).

The map $\tau$ is injective. Thus, if we had $\tau\left(  j\right)
=\tau\left(  i\right)  $, then we would have $j=i$, which would contradict
$j\neq i$. Hence, we cannot have $\tau\left(  j\right)  =\tau\left(  i\right)
$. Thus, we have $\tau\left(  j\right)  \neq\tau\left(  i\right)  $. In other
words, $\tau\left(  i\right)  \neq\tau\left(  j\right)  $. The same argument
(but applied to $\sigma$ instead of $\tau$) yields $\sigma\left(  i\right)
\neq\sigma\left(  j\right)  $.

Now, we have the logical implication%
\begin{equation}
\left(  \tau\left(  i\right)  >\tau\left(  j\right)  \right)
\ \Longrightarrow\ \left(  \sigma\left(  i\right)  >\sigma\left(  j\right)
\right)  . \label{pf.lem.sol.perm.lehmer.lili.c.c2.imp1}%
\end{equation}

[\textit{Proof of (\ref{pf.lem.sol.perm.lehmer.lili.c.c2.imp1}):} Assume that
$\tau\left(  i\right)  >\tau\left(  j\right)  $ holds. We must prove that
$\sigma\left(  i\right)  >\sigma\left(  j\right)  $ holds.

We have $\tau\left(  i\right)  >\tau\left(  j\right)  $, thus $\tau\left(
i\right)  \geq\tau\left(  j\right)  +1$ (since $\tau\left(  i\right)  $ and
$\tau\left(  j\right)  $ are integers). Hence, $\tau\left(  i\right)
-1\geq\tau\left(  j\right)  $.

We have $\tau\left(  i\right)  -1\leq\sigma\left(  i\right)  $, thus
$\sigma\left(  i\right)  \geq\tau\left(  i\right)  -1\geq\tau\left(  j\right)
\geq\sigma\left(  j\right)  $ (since $\sigma\left(  j\right)  \leq\tau\left(
j\right)  $). Combining this with $\sigma\left(  i\right)  \neq\sigma\left(
j\right)  $, we obtain $\sigma\left(  i\right)  >\sigma\left(  j\right)  $.

Now, forget our assumption that $\tau\left(  i\right)  >\tau\left(  j\right)
$. We thus have proven that if $\tau\left(  i\right)  >\tau\left(  j\right)  $
holds, then $\sigma\left(  i\right)  >\sigma\left(  j\right)  $ holds. Thus,
the implication (\ref{pf.lem.sol.perm.lehmer.lili.c.c2.imp1}) is proven.]

We also have the logical implication%
\begin{equation}
\left(  \sigma\left(  i\right)  >\sigma\left(  j\right)  \right)
\ \Longrightarrow\ \left(  \tau\left(  i\right)  >\tau\left(  j\right)
\right)  . \label{pf.lem.sol.perm.lehmer.lili.c.c2.imp2}%
\end{equation}

[\textit{Proof of (\ref{pf.lem.sol.perm.lehmer.lili.c.c2.imp2}):} Assume that
$\sigma\left(  i\right)  >\sigma\left(  j\right)  $ holds. We must prove that
$\tau\left(  i\right)  >\tau\left(  j\right)  $ holds.

We have $\sigma\left(  i\right)  >\sigma\left(  j\right)  $, thus
$\sigma\left(  i\right)  \geq\sigma\left(  j\right)  +1$ (since $\sigma\left(
i\right)  $ and $\sigma\left(  j\right)  $ are integers). Hence,
$\sigma\left(  i\right)  -1\geq\sigma\left(  j\right)  $.

We have $\tau\left(  j\right)  -1\leq\sigma\left(  j\right)  $, thus
$\sigma\left(  j\right)  \geq\tau\left(  j\right)  -1$ and therefore
$\sigma\left(  j\right)  +1\geq\tau\left(  j\right)  $.

But from $\sigma\left(  i\right)  \leq\tau\left(  i\right)  $, we obtain
$\tau\left(  i\right)  \geq\sigma\left(  i\right)  \geq\sigma\left(  j\right)
+1\geq\tau\left(  j\right)  $. Combining this with $\tau\left(  i\right)
\neq\tau\left(  j\right)  $, we obtain $\tau\left(  i\right)  >\tau\left(
j\right)  $.

Now, forget our assumption that $\sigma\left(  i\right)  >\sigma\left(
j\right)  $. We thus have proven that if $\sigma\left(  i\right)
>\sigma\left(  j\right)  $ holds, then $\tau\left(  i\right)  >\tau\left(
j\right)  $ holds. Thus, the implication
(\ref{pf.lem.sol.perm.lehmer.lili.c.c2.imp2}) is proven.]

Combining the two implications (\ref{pf.lem.sol.perm.lehmer.lili.c.c2.imp1})
and (\ref{pf.lem.sol.perm.lehmer.lili.c.c2.imp2}), we obtain the logical
equivalence $\left(  \tau\left(  i\right)  >\tau\left(  j\right)  \right)
\ \Longleftrightarrow\ \left(  \sigma\left(  i\right)  >\sigma\left(
j\right)  \right)  $. This proves Claim 2.]

Now, recall that $\ell_{i}\left(  \sigma\right)  $ is the number of all
$j\in\left\{  i+1,i+2,\ldots,n\right\}  $ such that $\sigma\left(  i\right)
>\sigma\left(  j\right)  $ (by the definition of $\ell_{i}\left(
\sigma\right)  $). In other words,%
\begin{equation}
\ell_{i}\left(  \sigma\right)  =\left\vert \left\{  j\in\left\{
i+1,i+2,\ldots,n\right\}  \ \mid\ \sigma\left(  i\right)  >\sigma\left(
j\right)  \right\}  \right\vert . \label{pf.lem.sol.perm.lehmer.lili.c.finsi}%
\end{equation}
The same argument (applied to $\tau$ instead of $\sigma$) yields%
\begin{equation}
\ell_{i}\left(  \tau\right)  =\left\vert \left\{  j\in\left\{  i+1,i+2,\ldots
,n\right\}  \ \mid\ \tau\left(  i\right)  >\tau\left(  j\right)  \right\}
\right\vert . \label{pf.lem.sol.perm.lehmer.lili.c.fintau}%
\end{equation}
But%
\begin{align*}
&  \left\{  j\in\left\{  i+1,i+2,\ldots,n\right\}  \ \mid\ \underbrace{\tau
\left(  i\right)  >\tau\left(  j\right)  }_{\substack{\Longleftrightarrow
\ \left(  \sigma\left(  i\right)  >\sigma\left(  j\right)  \right)
\\\text{(by Claim 2)}}}\right\} \\
&  =\left\{  j\in\left\{  i+1,i+2,\ldots,n\right\}  \ \mid\ \sigma\left(
i\right)  >\sigma\left(  j\right)  \right\}  .
\end{align*}
Hence, (\ref{pf.lem.sol.perm.lehmer.lili.c.fintau}) becomes%
\begin{align*}
\ell_{i}\left(  \tau\right)   &  =\left\vert \underbrace{\left\{  j\in\left\{
i+1,i+2,\ldots,n\right\}  \ \mid\ \tau\left(  i\right)  >\tau\left(  j\right)
\right\}  }_{=\left\{  j\in\left\{  i+1,i+2,\ldots,n\right\}  \ \mid
\ \sigma\left(  i\right)  >\sigma\left(  j\right)  \right\}  }\right\vert \\
&  =\left\vert \left\{  j\in\left\{  i+1,i+2,\ldots,n\right\}  \ \mid
\ \sigma\left(  i\right)  >\sigma\left(  j\right)  \right\}  \right\vert
=\ell_{i}\left(  \sigma\right)
\end{align*}
(by (\ref{pf.lem.sol.perm.lehmer.lili.c.finsi})). This proves Lemma
\ref{lem.sol.perm.lehmer.lili} \textbf{(b)}.

\textbf{(c)} Assume that we have (\ref{eq.lem.sol.perm.lehmer.lili.c.ass}). We
must prove that we have (\ref{eq.lem.sol.perm.lehmer.lili.c.claim}). In other
words, we must prove that $\tau\left(  q\right)  =q$ for each $q\in\left\{
1,2,\ldots,u\right\}  $.

So let $q\in\left\{  1,2,\ldots,u\right\}  $. We then must prove that
$\tau\left(  q\right)  =q$.

We are in one of the following two cases:

\textit{Case 1:} We have $q=u$.

\textit{Case 2:} We have $q\neq u$.

Let us first consider Case 1. In this case, we have $q=u$. Thus,
\begin{align*}
\tau\left(  \underbrace{q}_{=u}\right)   &  =\tau\left(  u\right)
=u\ \ \ \ \ \ \ \ \ \ \left(  \text{by Lemma \ref{lem.sol.perm.lehmer.lili}
\textbf{(a)}}\right) \\
&  =q.
\end{align*}
Hence, $\tau\left(  q\right)  =q$ is proven in Case 1.

Let us now consider Case 2. In this case, we have $q\neq u$. Combining
$q\in\left\{  1,2,\ldots,u\right\}  $ with $q\neq u$, we obtain $q\in\left\{
1,2,\ldots,u\right\}  \setminus\left\{  u\right\}  =\left\{  1,2,\ldots
,u-1\right\}  $. Thus, (\ref{eq.lem.sol.perm.lehmer.lili.c.ass}) yields
$\sigma\left(  q\right)  =q$. Also, from $q\in\left\{  1,2,\ldots,u-1\right\}
$, we obtain $q\leq u-1<u$.

Also,
\begin{align*}
q  &  \in\left\{  1,2,\ldots,u\right\}  \subseteq\left\{  1,2,\ldots
,n\right\}  \ \ \ \ \ \ \ \ \ \ \left(  \text{since }u\leq n\text{ (because
}u\in\left[  n\right]  =\left\{  1,2,\ldots,n\right\}  \text{)}\right) \\
&  =\left[  n\right]  .
\end{align*}
If we had $q\in\left\{  u,u+1,\ldots,v\right\}  $, then we would have $q\geq
u$, which would contradict $q<u$. Hence, we cannot have $q\in\left\{
u,u+1,\ldots,v\right\}  $. In other words, we have $q\notin\left\{
u,u+1,\ldots,v\right\}  $.

Combining $q\in\left[  n\right]  $ with $q\notin\left\{  u,u+1,\ldots
,v\right\}  $, we obtain $q\in\left[  n\right]  \setminus\left\{
u,u+1,\ldots,v\right\}  $. Hence, Lemma \ref{lem.sol.perm.lehmer.cuv.values}
\textbf{(a)} yields $c_{u,v}\left(  q\right)  =q$.

Now,%
\[
\underbrace{\tau}_{=\left(  c_{u,v}\right)  ^{-1}\circ\sigma}\left(  q\right)
=\left(  \left(  c_{u,v}\right)  ^{-1}\circ\sigma\right)  \left(  q\right)
=\left(  c_{u,v}\right)  ^{-1}\left(  \underbrace{\sigma\left(  q\right)
}_{=q}\right)  =\left(  c_{u,v}\right)  ^{-1}\left(  q\right)  =q
\]
(since $c_{u,v}\left(  q\right)  =q$). Hence, $\tau\left(  q\right)  =q$ is
proven in Case 2.

\begin{vershort}
We have now proven $\tau\left(  q\right)  =q$ in each of the two Cases 1 and
2. We thus conclude that $\tau\left(  q\right)  =q$ always holds.
\end{vershort}

\begin{verlong}
We have now proven $\tau\left(  q\right)  =q$ in each of the two Cases 1 and
2. Since these two Cases cover all possibilities, we thus conclude that
$\tau\left(  q\right)  =q$ always holds.
\end{verlong}

Thus, we have shown that we have (\ref{eq.lem.sol.perm.lehmer.lili.c.claim}).
This proves Lemma \ref{lem.sol.perm.lehmer.lili} \textbf{(c)}.
\end{proof}

\begin{lemma}
\label{lem.sol.perm.lehmer.rothe.liw0}Let $n\in\mathbb{N}$. Define the
permutation $w_{0}\in S_{n}$ as in Exercise \ref{exe.ps2.2.4} \textbf{(c)}.
Then,%
\begin{equation}
\ell_{i}\left(  w_{0}\right)  =n-i\ \ \ \ \ \ \ \ \ \ \text{for each }%
i\in\left\{  1,2,\ldots,n\right\}  . \label{sol.perm.lehmer.rothe.e.liw0}%
\end{equation}

\end{lemma}

\begin{proof}
[Proof of Lemma \ref{lem.sol.perm.lehmer.rothe.liw0}.]Let $i\in\left\{
1,2,\ldots,n\right\}  $. Thus, $i\geq1$ and $i\leq n$ and $i\in\left\{
1,2,\ldots,n\right\}  =\left[  n\right]  $ (since $\left[  n\right]  =\left\{
1,2,\ldots,n\right\}  $ (by the definition of $\left[  n\right]  $)).

We know that $\ell_{i}\left(  w_{0}\right)  $ is the number of all
$j\in\left\{  i+1,i+2,\ldots,n\right\}  $ such that $w_{0}\left(  i\right)
>w_{0}\left(  j\right)  $ (by the definition of $\ell_{i}\left(  w_{0}\right)
$). In other words,%
\begin{equation}
\ell_{i}\left(  w_{0}\right)  =\left\vert \left\{  j\in\left\{  i+1,i+2,\ldots
,n\right\}  \ \mid\ w_{0}\left(  i\right)  >w_{0}\left(  j\right)  \right\}
\right\vert . \label{sol.perm.lehmer.rothe.e.liw0.pf.1}%
\end{equation}

\begin{vershort}
Now, let $k\in\left\{  i+1,i+2,\ldots,n\right\}  $. Thus, $k\geq i+1>i$. The
definition of $w_{0}$ yields $w_{0}\left(  i\right)  =n+1-i$ and $w_{0}\left(
k\right)  =n+1-\underbrace{k}_{>i}<n+1-i=w_{0}\left(  i\right)  $, so that
$w_{0}\left(  i\right)  >w_{0}\left(  k\right)  $. Hence, $k$ is a
$j\in\left\{  i+1,i+2,\ldots,n\right\}  $ satisfying $w_{0}\left(  i\right)
>w_{0}\left(  j\right)  $. In other words, $k\in\left\{  j\in\left\{
i+1,i+2,\ldots,n\right\}  \ \mid\ w_{0}\left(  i\right)  >w_{0}\left(
j\right)  \right\}  $.
\end{vershort}

\begin{verlong}
Now, let $k\in\left\{  i+1,i+2,\ldots,n\right\}  $. Thus, $k\geq i+1>i$. Also,
$k\in\left\{  i+1,i+2,\ldots,n\right\}  \subseteq\left\{  1,2,\ldots
,n\right\}  $ (since $i+1\geq i\geq1$). But the definition of $w_{0}$ yields
$w_{0}\left(  i\right)  =n+1-i$ and $w_{0}\left(  k\right)  =n+1-k$. Hence,
$w_{0}\left(  k\right)  =n+1-\underbrace{k}_{>i}<n+1-i=w_{0}\left(  i\right)
$, so that $w_{0}\left(  i\right)  >w_{0}\left(  k\right)  $. Hence, $k$ is an
element of $\left\{  i+1,i+2,\ldots,n\right\}  $ and satisfies $w_{0}\left(
i\right)  >w_{0}\left(  k\right)  $. In other words, $k$ is a $j\in\left\{
i+1,i+2,\ldots,n\right\}  $ satisfying $w_{0}\left(  i\right)  >w_{0}\left(
j\right)  $. In other words, $k\in\left\{  j\in\left\{  i+1,i+2,\ldots
,n\right\}  \ \mid\ w_{0}\left(  i\right)  >w_{0}\left(  j\right)  \right\}  $.
\end{verlong}

Now, forget that we fixed $k$. We thus have proven that \newline$k\in\left\{
j\in\left\{  i+1,i+2,\ldots,n\right\}  \ \mid\ w_{0}\left(  i\right)
>w_{0}\left(  j\right)  \right\}  $ for each $k\in\left\{  i+1,i+2,\ldots
,n\right\}  $. In other words, we have%
\[
\left\{  i+1,i+2,\ldots,n\right\}  \subseteq\left\{  j\in\left\{
i+1,i+2,\ldots,n\right\}  \ \mid\ w_{0}\left(  i\right)  >w_{0}\left(
j\right)  \right\}  .
\]
Combining this inclusion with the (obvious) inclusion%
\[
\left\{  j\in\left\{  i+1,i+2,\ldots,n\right\}  \ \mid\ w_{0}\left(  i\right)
>w_{0}\left(  j\right)  \right\}  \subseteq\left\{  i+1,i+2,\ldots,n\right\}
,
\]
we obtain%
\[
\left\{  j\in\left\{  i+1,i+2,\ldots,n\right\}  \ \mid\ w_{0}\left(  i\right)
>w_{0}\left(  j\right)  \right\}  =\left\{  i+1,i+2,\ldots,n\right\}  .
\]
Hence, (\ref{sol.perm.lehmer.rothe.e.liw0.pf.1}) becomes%
\begin{align*}
\ell_{i}\left(  w_{0}\right)   &  =\left\vert \underbrace{\left\{
j\in\left\{  i+1,i+2,\ldots,n\right\}  \ \mid\ w_{0}\left(  i\right)
>w_{0}\left(  j\right)  \right\}  }_{=\left\{  i+1,i+2,\ldots,n\right\}
}\right\vert =\left\vert \left\{  i+1,i+2,\ldots,n\right\}  \right\vert \\
&  =n-i\ \ \ \ \ \ \ \ \ \ \left(  \text{since }i\leq n\right)  .
\end{align*}
This proves Lemma \ref{lem.sol.perm.lehmer.rothe.liw0}.
\end{proof}

Next, let us state an induction principle that will come useful in our
solution of Exercise \ref{exe.perm.lehmer.rothe} \textbf{(c)}:

\begin{theorem}
\label{thm.sol.perm.lehmer.IPgh-1}Let $g\in\mathbb{Z}$ and $h\in\mathbb{Z}$.
For each $p\in\left\{  g,g+1,\ldots,h\right\}  $, let $\mathcal{A}\left(
p\right)  $ be a logical statement.

Assume the following:

\begin{statement}
\textit{Assumption 1:} If $g\leq h$, then the statement $\mathcal{A}\left(
h\right)  $ holds.
\end{statement}

\begin{statement}
\textit{Assumption 2:} If $u\in\left\{  g+1,g+2,\ldots,h\right\}  $ is such
that $\mathcal{A}\left(  u\right)  $ holds, then $\mathcal{A}\left(
u-1\right)  $ also holds.
\end{statement}

Then, $\mathcal{A}\left(  p\right)  $ holds for each $p\in\left\{
g,g+1,\ldots,h\right\}  $.
\end{theorem}

\begin{proof}
[Proof of Theorem \ref{thm.sol.perm.lehmer.IPgh-1}.]Theorem
\ref{thm.sol.perm.lehmer.IPgh-1} is exactly Theorem \ref{thm.ind.IPgh-},
except that some names have been changed:

\begin{itemize}
\item The variable $n$ has been renamed as $p$.

\item The variable $m$ has been renamed as $u$.
\end{itemize}

Thus, Theorem \ref{thm.sol.perm.lehmer.IPgh-1} holds (since Theorem
\ref{thm.ind.IPgh-} holds).
\end{proof}

We shall only need the particular case of Theorem
\ref{thm.sol.perm.lehmer.IPgh-1} in which $g=0$ and $h=n$:

\begin{verlong}
\begin{corollary}
\label{cor.sol.perm.lehmer.IP0n-1}Let $n\in\mathbb{Z}$. For each $p\in\left\{
0,0+1,\ldots,n\right\}  $, let $\mathcal{A}\left(  p\right)  $ be a logical statement.

Assume the following:

\begin{statement}
\textit{Assumption 1:} If $0\leq n$, then the statement $\mathcal{A}\left(
n\right)  $ holds.
\end{statement}

\begin{statement}
\textit{Assumption 2:} If $u\in\left\{  0+1,0+2,\ldots,n\right\}  $ is such
that $\mathcal{A}\left(  u\right)  $ holds, then $\mathcal{A}\left(
u-1\right)  $ also holds.
\end{statement}

Then, $\mathcal{A}\left(  p\right)  $ holds for each $p\in\left\{
0,0+1,\ldots,n\right\}  $.
\end{corollary}

\begin{proof}
[Proof of Corollary \ref{cor.sol.perm.lehmer.IP0n-1}.]Corollary
\ref{cor.sol.perm.lehmer.IP0n-1} is just the particular case of Theorem
\ref{thm.sol.perm.lehmer.IPgh-1} obtained by taking $g=0$ and $h=n$.
\end{proof}

Let us restate Corollary \ref{cor.sol.perm.lehmer.IP0n-1} with slightly
friendlier notations:
\end{verlong}

\begin{corollary}
\label{cor.sol.perm.lehmer.IP0n-1f}Let $n\in\mathbb{Z}$. For each
$p\in\left\{  0,1,\ldots,n\right\}  $, let $\mathcal{A}\left(  p\right)  $ be
a logical statement.

Assume the following:

\begin{statement}
\textit{Assumption 1:} If $0\leq n$, then the statement $\mathcal{A}\left(
n\right)  $ holds.
\end{statement}

\begin{statement}
\textit{Assumption 2:} If $u\in\left\{  1,2,\ldots,n\right\}  $ is such that
$\mathcal{A}\left(  u\right)  $ holds, then $\mathcal{A}\left(  u-1\right)  $
also holds.
\end{statement}

Then, $\mathcal{A}\left(  p\right)  $ holds for each $p\in\left\{
0,1,\ldots,n\right\}  $.
\end{corollary}

\begin{vershort}
\begin{proof}
[Proof of Corollary \ref{cor.sol.perm.lehmer.IP0n-1f}.]Corollary
\ref{cor.sol.perm.lehmer.IP0n-1f} is just the particular case of Theorem
\ref{thm.sol.perm.lehmer.IPgh-1} obtained by taking $g=0$ and $h=n$.
\end{proof}
\end{vershort}

\begin{verlong}
\begin{proof}
[Proof of Corollary \ref{cor.sol.perm.lehmer.IP0n-1f}.]Clearly,%
\[
\left\{  0,0+1,\ldots,n\right\}  =\left\{  0,1,\ldots,n\right\}
\ \ \ \ \ \ \ \ \ \ \text{and}\ \ \ \ \ \ \ \ \ \ \left\{  0+1,0+2,\ldots
,n\right\}  =\left\{  1,2,\ldots,n\right\}  .
\]
Rewriting Corollary \ref{cor.sol.perm.lehmer.IP0n-1} using these two
equalities, we obtain precisely Corollary \ref{cor.sol.perm.lehmer.IP0n-1f}.
Thus, Corollary \ref{cor.sol.perm.lehmer.IP0n-1f} is proven.
\end{proof}
\end{verlong}

At this point we have proven enough auxiliary facts to render the solution of
Exercise \ref{exe.perm.lehmer.rothe} fairly anodyne:

\begin{proof}
[Solution to Exercise \ref{exe.perm.lehmer.rothe}.]Let us first observe that
the definition of $c_{u,v}$ given in Exercise \ref{exe.perm.lehmer.rothe} is
exactly the definition of $c_{u,v}$ given in Definition
\ref{def.sol.perm.c=ttt.cuv}. Thus, there is no conflict of notation.

\begin{verlong}
Recall that $\left[  n\right]  =\left\{  1,2,\ldots,n\right\}  $ (by the
definition of $\left[  n\right]  $).
\end{verlong}

\textbf{(a)} Let $i\in\left[  n\right]  $. Lemma \ref{lem.sol.perm.lehmer.ab}
\textbf{(a)} yields that the permutation $c_{i,i+\ell_{i}\left(
\sigma\right)  }\in S_{n}$ is well-defined. Thus, the permutation $a_{i}\in
S_{n}$ is well-defined (since $a_{i}$ was defined by $a_{i}=c_{i,i+\ell
_{i}\left(  \sigma\right)  }$). This solves Exercise
\ref{exe.perm.lehmer.rothe} \textbf{(a)}.

\textbf{(b)} Let $i\in\left[  n\right]  $. Set $i^{\prime}=i+\ell_{i}\left(
\sigma\right)  $. We must show that $a_{i}=s_{i^{\prime}-1}\circ s_{i^{\prime
}-2}\circ\cdots\circ s_{i}$.

The definition of $a_{i}$ yields $a_{i}=c_{i,i+\ell_{i}\left(  \sigma\right)
}=s_{i^{\prime}-1}\circ s_{i^{\prime}-2}\circ\cdots\circ s_{i}$ (by Lemma
\ref{lem.sol.perm.lehmer.ab} \textbf{(b)}). This solves Exercise
\ref{exe.perm.lehmer.rothe} \textbf{(b)}.

\textbf{(c)} Forget that we fixed $\sigma$ (but let us leave $n$ fixed).

Let us define some further notations. If $\sigma\in S_{n}$ and $i\in\left[
n\right]  $ are arbitrary, then we define a permutation $a_{i}^{\left(
\sigma\right)  }\in S_{n}$ by%
\begin{equation}
a_{i}^{\left(  \sigma\right)  }=c_{i,i+\ell_{i}\left(  \sigma\right)  }.
\label{sol.perm.lehmer.rothe.c.aisig=}%
\end{equation}
(This is well-defined, because Lemma \ref{lem.sol.perm.lehmer.ab} \textbf{(a)}
yields that the permutation $c_{i,i+\ell_{i}\left(  \sigma\right)  }\in S_{n}$
is well-defined.) Thus, $a_{i}^{\left(  \sigma\right)  }$ is precisely what
was called $a_{i}$ in Exercise \ref{exe.perm.lehmer.rothe}, but we have
renamed it $a_{i}^{\left(  \sigma\right)  }$ in order to make its dependence
on $\sigma$ explicit. (This allows us to simultaneously consider
$a_{i}^{\left(  \sigma\right)  }$ and $a_{i}^{\left(  \tau\right)  }$ for two
different permutations $\sigma$ and $\tau$.)

If $p\in\left\{  0,1,\ldots,n\right\}  $ and $\sigma\in S_{n}$, then we say
that $\sigma$ is $p$\textit{-lazy} if we have%
\[
\left(  \sigma\left(  q\right)  =q\text{ for each }q\in\left\{  1,2,\ldots
,p\right\}  \right)  .
\]
Note that
\begin{equation}
\text{every permutation }\sigma\in S_{n}\text{ is }0\text{-lazy}
\label{sol.perm.lehmer.rothe.c.0-lazy}%
\end{equation}
\footnote{\textit{Proof of (\ref{sol.perm.lehmer.rothe.c.0-lazy}):} Let
$\sigma\in S_{n}$. We must prove that $\sigma$ is $0$-lazy.
\par
The statement $\left(  \sigma\left(  q\right)  =q\text{ for each }q\in\left\{
1,2,\ldots,0\right\}  \right)  $ is vacuously true (since there exists no
$q\in\left\{  1,2,\ldots,0\right\}  $ (because $\left\{  1,2,\ldots,0\right\}
=\varnothing$)), and therefore is true. In other words, $\sigma$ is $0$-lazy
(because $\sigma$ is $0$-lazy if and only if we have $\left(  \sigma\left(
q\right)  =q\text{ for each }q\in\left\{  1,2,\ldots,0\right\}  \right)  $ (by
the definition of \textquotedblleft$0$-lazy\textquotedblright)). This proves
(\ref{sol.perm.lehmer.rothe.c.0-lazy}).}. On the other hand,%
\begin{equation}
\text{if a permutation }\sigma\in S_{n}\text{ is }n\text{-lazy, then }%
\sigma=\operatorname*{id} \label{sol.perm.lehmer.rothe.c.n-lazy}%
\end{equation}
\footnote{\textit{Proof of (\ref{sol.perm.lehmer.rothe.c.n-lazy}):} Let
$\sigma\in S_{n}$ be a permutation that is $n$-lazy. We must prove that
$\sigma=\operatorname*{id}$.
\par
We have $\sigma\in S_{n}$. In other words, $\sigma$ is a permutation of
$\left\{  1,2,\ldots,n\right\}  $ (since $S_{n}$ is the set of all
permutations of $\left\{  1,2,\ldots,n\right\}  $). In other words, $\sigma$
is a bijection from $\left\{  1,2,\ldots,n\right\}  $ to $\left\{
1,2,\ldots,n\right\}  $.
\par
But $\sigma$ is $n$-lazy if and only if we have $\left(  \sigma\left(
q\right)  =q\text{ for each }q\in\left\{  1,2,\ldots,n\right\}  \right)  $ (by
the definition of \textquotedblleft$n$-lazy\textquotedblright). Hence, we have
$\left(  \sigma\left(  q\right)  =q\text{ for each }q\in\left\{
1,2,\ldots,n\right\}  \right)  $ (since $\sigma$ is $n$-lazy). Thus, for each
$q\in\left\{  1,2,\ldots,n\right\}  $, we have $\sigma\left(  q\right)
=q=\operatorname*{id}\left(  q\right)  $. Hence, $\sigma=\operatorname*{id}$
(since both $\sigma$ and $\operatorname*{id}$ are maps from $\left\{
1,2,\ldots,n\right\}  $ to $\left\{  1,2,\ldots,n\right\}  $). This proves
(\ref{sol.perm.lehmer.rothe.c.n-lazy}).}.

For each $p\in\left\{  0,1,\ldots,n\right\}  $, let us define a logical
statement $\mathcal{A}\left(  p\right)  $ as follows:

\begin{statement}
\textit{Statement $\mathcal{A}\left(  p\right)  $:} Every $p$-lazy permutation
$\sigma\in S_{n}$ satisfies%
\[
\sigma=a_{p+1}^{\left(  \sigma\right)  }\circ a_{p+2}^{\left(  \sigma\right)
}\circ\cdots\circ a_{n}^{\left(  \sigma\right)  }.
\]

\end{statement}

Our goal is to prove that $\mathcal{A}\left(  p\right)  $ holds for each
$p\in\left\{  0,1,\ldots,n\right\}  $. Once this is done, we will immediately
conclude that $\mathcal{A}\left(  0\right)  $ holds, which means that Exercise
\ref{exe.perm.lehmer.rothe} \textbf{(c)} holds (since the claim of Exercise
\ref{exe.perm.lehmer.rothe} \textbf{(c)} is more or less obviously equivalent
to $\mathcal{A}\left(  0\right)  $).

In order to prove that $\mathcal{A}\left(  p\right)  $ holds for each
$p\in\left\{  0,1,\ldots,n\right\}  $, we shall use Corollary
\ref{cor.sol.perm.lehmer.IP0n-1f}. In order to do so, we need to prove that
Assumptions 1 and 2 of Corollary \ref{cor.sol.perm.lehmer.IP0n-1f} hold for
these statements $\mathcal{A}\left(  0\right)  ,\mathcal{A}\left(  1\right)
,\ldots,\mathcal{A}\left(  n\right)  $. Let us therefore prove that these two
assumptions hold:

[\textit{Proof of Assumption 1:} Assume that $0\leq n$. (This is clearly true
anyway, since $n\in\mathbb{N}$.)

Let $\sigma\in S_{n}$ be an $n$-lazy permutation. Thus, $\sigma
=\operatorname*{id}$ (by (\ref{sol.perm.lehmer.rothe.c.n-lazy})). Comparing
this with%
\[
a_{n+1}^{\left(  \sigma\right)  }\circ a_{n+2}^{\left(  \sigma\right)  }%
\circ\cdots\circ a_{n}^{\left(  \sigma\right)  }=\left(  \text{an empty
composition of permutations in }S_{n}\right)  =\operatorname*{id},
\]
we obtain $\sigma=a_{n+1}^{\left(  \sigma\right)  }\circ a_{n+2}^{\left(
\sigma\right)  }\circ\cdots\circ a_{n}^{\left(  \sigma\right)  }$.

Now, forget that we fixed $\sigma$. We thus have shown that every $n$-lazy
permutation $\sigma\in S_{n}$ satisfies $\sigma=a_{n+1}^{\left(
\sigma\right)  }\circ a_{n+2}^{\left(  \sigma\right)  }\circ\cdots\circ
a_{n}^{\left(  \sigma\right)  }$. But this is precisely the statement
$\mathcal{A}\left(  n\right)  $. Hence, we have shown that the statement
$\mathcal{A}\left(  n\right)  $ holds. This concludes the proof of Assumption 1.]

[\textit{Proof of Assumption 2:} Let $u\in\left\{  1,2,\ldots,n\right\}  $ be
such that $\mathcal{A}\left(  u\right)  $ holds. We must prove that
$\mathcal{A}\left(  u-1\right)  $ also holds.

We have $u\in\left\{  1,2,\ldots,n\right\}  =\left[  n\right]  $.

We have assumed that $\mathcal{A}\left(  u\right)  $ holds. In other words,
\begin{equation}
\text{every }u\text{-lazy permutation }\sigma\in S_{n}\text{ satisfies }%
\sigma=a_{u+1}^{\left(  \sigma\right)  }\circ a_{u+2}^{\left(  \sigma\right)
}\circ\cdots\circ a_{n}^{\left(  \sigma\right)  }
\label{sol.perm.lehmer.rothe.c.a2.pf.IH}%
\end{equation}
(because this is what the statement $\mathcal{A}\left(  u\right)  $ says).

Now, we will try to prove that $\mathcal{A}\left(  u-1\right)  $ holds.

Let $\sigma\in S_{n}$ be a $\left(  u-1\right)  $-lazy permutation. Thus,
$\sigma$ is $\left(  u-1\right)  $-lazy. In other words,%
\begin{equation}
\left(  \sigma\left(  q\right)  =q\text{ for each }q\in\left\{  1,2,\ldots
,u-1\right\}  \right)  \label{sol.perm.lehmer.rothe.c.a2.pf.sigq}%
\end{equation}
\footnote{because $\sigma$ is $\left(  u-1\right)  $-lazy if and only if
$\left(  \sigma\left(  q\right)  =q\text{ for each }q\in\left\{
1,2,\ldots,u-1\right\}  \right)  $ (by the definition of \textquotedblleft%
$\left(  u-1\right)  $-lazy\textquotedblright)}. Hence, Lemma
\ref{lem.sol.perm.lehmer.lu-if-fix} yields $u+\ell_{u}\left(  \sigma\right)
=\sigma\left(  u\right)  $.

\begin{verlong}
Recall that $S_{n}$ is the set of all permutations of the set $\left\{
1,2,\ldots,n\right\}  $. In other words, $S_{n}$ is the set of all
permutations of the set $\left[  n\right]  $ (since $\left[  n\right]
=\left\{  1,2,\ldots,n\right\}  $). We have $\sigma\in S_{n}$. In other words,
$\sigma$ is a permutation of $\left[  n\right]  $ (since $S_{n}$ is the set of
all permutations of the set $\left[  n\right]  $). In other words, $\sigma$ is
a bijective map $\left[  n\right]  \rightarrow\left[  n\right]  $. Hence,
$\sigma\left(  u\right)  \in\left[  n\right]  $.
\end{verlong}

Define $v\in\left[  n\right]  $ by $v=\sigma\left(  u\right)  $. (This is
well-defined, since $\sigma\left(  u\right)  \in\left[  n\right]  $.) Then,
$u+\ell_{u}\left(  \sigma\right)  =\sigma\left(  u\right)  =v$ (since
$v=\sigma\left(  u\right)  $).

\begin{verlong}
Note that $\ell_{u}\left(  \sigma\right)  $ is the number of all $j\in\left\{
u+1,u+2,\ldots,n\right\}  $ such that $\sigma\left(  u\right)  >\sigma\left(
j\right)  $ (by the definition of $\ell_{u}\left(  \sigma\right)  $). Thus,
$\ell_{u}\left(  \sigma\right)  $ is a nonnegative integer; hence, $\ell
_{u}\left(  \sigma\right)  \geq0$.
\end{verlong}

From $u+\ell_{u}\left(  \sigma\right)  =v$, we obtain $v=u+\underbrace{\ell
_{u}\left(  \sigma\right)  }_{\geq0}\geq u$, so that $u\leq v$.

Now, the definition of $a_{u}^{\left(  \sigma\right)  }$ yields%
\begin{equation}
a_{u}^{\left(  \sigma\right)  }=c_{u,u+\ell_{u}\left(  \sigma\right)
}=c_{u,v}\ \ \ \ \ \ \ \ \ \ \left(  \text{since }u+\ell_{u}\left(
\sigma\right)  =v\right)  . \label{sol.perm.lehmer.rothe.c.a2.pf.ausi=}%
\end{equation}

Define a permutation $\tau\in S_{n}$ by $\tau=\left(  c_{u,v}\right)
^{-1}\circ\sigma$. Hence, Lemma \ref{lem.sol.perm.lehmer.lili} \textbf{(b)}
yields that%
\begin{equation}
\ell_{i}\left(  \tau\right)  =\ell_{i}\left(  \sigma\right)
\ \ \ \ \ \ \ \ \ \ \text{for each }i\in\left\{  u+1,u+2,\ldots,n\right\}  .
\label{sol.perm.lehmer.rothe.c.a2.pf.li=li}%
\end{equation}
Thus, each $i\in\left\{  u+1,u+2,\ldots,n\right\}  $ satisfies%
\begin{align*}
a_{i}^{\left(  \tau\right)  }  &  =c_{i,i+\ell_{i}\left(  \tau\right)
}\ \ \ \ \ \ \ \ \ \ \left(  \text{by the definition of }a_{i}^{\left(
\tau\right)  }\right) \\
&  =c_{i,i+\ell_{i}\left(  \sigma\right)  }\ \ \ \ \ \ \ \ \ \ \left(
\text{since }\ell_{i}\left(  \tau\right)  =\ell_{i}\left(  \sigma\right)
\text{ (by (\ref{sol.perm.lehmer.rothe.c.a2.pf.li=li}))}\right) \\
&  =a_{i}^{\left(  \sigma\right)  }\ \ \ \ \ \ \ \ \ \ \left(  \text{by
(\ref{sol.perm.lehmer.rothe.c.aisig=})}\right)  .
\end{align*}
In other words, $\left(  a_{u+1}^{\left(  \tau\right)  },a_{u+2}^{\left(
\tau\right)  },\ldots,a_{n}^{\left(  \tau\right)  }\right)  =\left(
a_{u+1}^{\left(  \sigma\right)  },a_{u+2}^{\left(  \sigma\right)  }%
,\ldots,a_{n}^{\left(  \sigma\right)  }\right)  $. Hence,%
\begin{equation}
a_{u+1}^{\left(  \tau\right)  }\circ a_{u+2}^{\left(  \tau\right)  }%
\circ\cdots\circ a_{n}^{\left(  \tau\right)  }=a_{u+1}^{\left(  \sigma\right)
}\circ a_{u+2}^{\left(  \sigma\right)  }\circ\cdots\circ a_{n}^{\left(
\sigma\right)  }. \label{sol.perm.lehmer.rothe.c.a2.pf.a=a}%
\end{equation}

Recall that $\left(  \sigma\left(  q\right)  =q\text{ for each }q\in\left\{
1,2,\ldots,u-1\right\}  \right)  $. Hence, Lemma
\ref{lem.sol.perm.lehmer.lili} \textbf{(c)} yields that $\left(  \tau\left(
q\right)  =q\text{ for each }q\in\left\{  1,2,\ldots,u\right\}  \right)  $. In
other words, the permutation $\tau$ is $u$-lazy\footnote{because $\tau$ is
$u$-lazy if and only if $\left(  \tau\left(  q\right)  =q\text{ for each }%
q\in\left\{  1,2,\ldots,u\right\}  \right)  $ (by the definition of
\textquotedblleft$u$-lazy\textquotedblright)}. Therefore,
(\ref{sol.perm.lehmer.rothe.c.a2.pf.IH}) (applied to $\tau$ instead of
$\sigma$) yields%
\begin{equation}
\tau=a_{u+1}^{\left(  \tau\right)  }\circ a_{u+2}^{\left(  \tau\right)  }%
\circ\cdots\circ a_{n}^{\left(  \tau\right)  }=a_{u+1}^{\left(  \sigma\right)
}\circ a_{u+2}^{\left(  \sigma\right)  }\circ\cdots\circ a_{n}^{\left(
\sigma\right)  } \label{sol.perm.lehmer.rothe.c.a2.pf.tau=aaa}%
\end{equation}
(by (\ref{sol.perm.lehmer.rothe.c.a2.pf.a=a})).

Now, $c_{u,v}\circ\underbrace{\tau}_{=\left(  c_{u,v}\right)  ^{-1}\circ
\sigma}=\underbrace{c_{u,v}\circ\left(  c_{u,v}\right)  ^{-1}}%
_{=\operatorname*{id}}\circ\sigma=\sigma$, so that%
\begin{align*}
\sigma &  =\underbrace{c_{u,v}}_{\substack{=a_{u}^{\left(  \sigma\right)
}\\\text{(by (\ref{sol.perm.lehmer.rothe.c.a2.pf.ausi=}))}}}\circ
\underbrace{\tau}_{\substack{=a_{u+1}^{\left(  \sigma\right)  }\circ
a_{u+2}^{\left(  \sigma\right)  }\circ\cdots\circ a_{n}^{\left(
\sigma\right)  }\\\text{(by (\ref{sol.perm.lehmer.rothe.c.a2.pf.tau=aaa}))}%
}}=a_{u}^{\left(  \sigma\right)  }\circ\left(  a_{u+1}^{\left(  \sigma\right)
}\circ a_{u+2}^{\left(  \sigma\right)  }\circ\cdots\circ a_{n}^{\left(
\sigma\right)  }\right) \\
&  =a_{u}^{\left(  \sigma\right)  }\circ a_{u+1}^{\left(  \sigma\right)
}\circ\cdots\circ a_{n}^{\left(  \sigma\right)  }=a_{\left(  u-1\right)
+1}^{\left(  \sigma\right)  }\circ a_{\left(  u-1\right)  +2}^{\left(
\sigma\right)  }\circ\cdots\circ a_{n}^{\left(  \sigma\right)  }%
\end{align*}
(since $u=\left(  u-1\right)  +1$ and $u+1=\left(  u-1\right)  +2$).

Now, forget that we fixed $\sigma$. We thus have proven that
\[
\text{every }\left(  u-1\right)  \text{-lazy permutation }\sigma\in
S_{n}\text{ satisfies }\sigma=a_{\left(  u-1\right)  +1}^{\left(
\sigma\right)  }\circ a_{\left(  u-1\right)  +2}^{\left(  \sigma\right)
}\circ\cdots\circ a_{n}^{\left(  \sigma\right)  }.
\]
But this is exactly the statement $\mathcal{A}\left(  u-1\right)  $. Thus, we
have proven that the statement $\mathcal{A}\left(  u-1\right)  $ holds. This
proves Assumption 2.]

We have now verified that both Assumptions 1 and 2 of Corollary
\ref{cor.sol.perm.lehmer.IP0n-1f} hold. Hence, Corollary
\ref{cor.sol.perm.lehmer.IP0n-1f} shows that $\mathcal{A}\left(  p\right)  $
holds for each $p\in\left\{  0,1,\ldots,n\right\}  $. Applying this to $p=0$,
we conclude that $\mathcal{A}\left(  0\right)  $ holds (since $0\in\left\{
0,1,\ldots,n\right\}  $). In other words,%
\begin{equation}
\text{every }0\text{-lazy permutation }\sigma\in S_{n}\text{ satisfies }%
\sigma=a_{0+1}^{\left(  \sigma\right)  }\circ a_{0+2}^{\left(  \sigma\right)
}\circ\cdots\circ a_{n}^{\left(  \sigma\right)  }
\label{sol.perm.lehmer.rothe.c.A0}%
\end{equation}
(because this is what the statement $\mathcal{A}\left(  0\right)  $ says).

Now, let $\sigma\in S_{n}$. Let us use the notation $a_{i}$ (for $i\in\left[
n\right]  $) as defined in Exercise \ref{exe.perm.lehmer.rothe}. The
permutation $\sigma\in S_{n}$ is $0$-lazy (by
(\ref{sol.perm.lehmer.rothe.c.0-lazy})). Thus,
(\ref{sol.perm.lehmer.rothe.c.A0}) yields%
\[
\sigma=a_{0+1}^{\left(  \sigma\right)  }\circ a_{0+2}^{\left(  \sigma\right)
}\circ\cdots\circ a_{n}^{\left(  \sigma\right)  }=a_{1}^{\left(
\sigma\right)  }\circ a_{2}^{\left(  \sigma\right)  }\circ\cdots\circ
a_{n}^{\left(  \sigma\right)  }.
\]
But each $i\in\left[  n\right]  $ satisfies%
\begin{align*}
a_{i}^{\left(  \sigma\right)  }  &  =c_{i,i+\ell_{i}\left(  \sigma\right)
}\ \ \ \ \ \ \ \ \ \ \left(  \text{by (\ref{sol.perm.lehmer.rothe.c.aisig=}%
)}\right) \\
&  =a_{i}\ \ \ \ \ \ \ \ \ \ \left(  \text{since }a_{i}=c_{i,i+\ell_{i}\left(
\sigma\right)  }\text{ (by the definition of }a_{i}\text{)}\right)  .
\end{align*}
In other words, each $i\in\left\{  1,2,\ldots,n\right\}  $ satisfies
$a_{i}^{\left(  \sigma\right)  }=a_{i}$ (since $\left[  n\right]  =\left\{
1,2,\ldots,n\right\}  $). In other words, $\left(  a_{1}^{\left(
\sigma\right)  },a_{2}^{\left(  \sigma\right)  },\ldots,a_{n}^{\left(
\sigma\right)  }\right)  =\left(  a_{1},a_{2},\ldots,a_{n}\right)  $. Hence,
$a_{1}^{\left(  \sigma\right)  }\circ a_{2}^{\left(  \sigma\right)  }%
\circ\cdots\circ a_{n}^{\left(  \sigma\right)  }=a_{1}\circ a_{2}\circ
\cdots\circ a_{n}$. Thus, $\sigma=a_{1}^{\left(  \sigma\right)  }\circ
a_{2}^{\left(  \sigma\right)  }\circ\cdots\circ a_{n}^{\left(  \sigma\right)
}=a_{1}\circ a_{2}\circ\cdots\circ a_{n}$. This solves Exercise
\ref{exe.perm.lehmer.rothe} \textbf{(c)}.

\textbf{(d)} \textit{Second solution to Exercise \ref{exe.ps2.2.5}
\textbf{(e)}:} In the following, a \textit{simple transposition} shall mean a
permutation of the form $s_{k}$ (with $k\in\left\{  1,2,\ldots,n-1\right\}
$). Thus, $s_{k}$ is a simple transposition for each $k\in\left\{
1,2,\ldots,n-1\right\}  $.

For each $i\in\left[  n\right]  $, we define a permutation $a_{i}\in S_{n}$ as
in Exercise \ref{exe.perm.lehmer.rothe}. Then, we claim that each
$i\in\left\{  1,2,\ldots,n\right\}  $ satisfies%
\begin{equation}
a_{i}=\left(  \text{a composition of }\ell_{i}\left(  \sigma\right)  \text{
simple transpositions}\right)  \label{sol.perm.lehmer.rothe.d.ai=}%
\end{equation}
(this equality is supposed to mean that $a_{i}$ is a composition of $i$ simple
transpositions; it does not mean that \textbf{every} composition of $i$ simple
transpositions is $a_{i}$).

[\textit{Proof of (\ref{sol.perm.lehmer.rothe.d.ai=}):} Let $i\in\left\{
1,2,\ldots,n\right\}  $. Thus, $i\in\left\{  1,2,\ldots,n\right\}  =\left[
n\right]  $. Set $i^{\prime}=i+\ell_{i}\left(  \sigma\right)  $. Then,
Exercise \ref{exe.perm.lehmer.rothe} \textbf{(b)} yields $a_{i}=s_{i^{\prime
}-1}\circ s_{i^{\prime}-2}\circ\cdots\circ s_{i}$.

\begin{verlong}
Note that $\ell_{i}\left(  \sigma\right)  $ is the number of all $j\in\left\{
i+1,i+2,\ldots,n\right\}  $ such that $\sigma\left(  i\right)  >\sigma\left(
j\right)  $ (by the definition of $\ell_{i}\left(  \sigma\right)  $). Thus,
$\ell_{i}\left(  \sigma\right)  $ is a nonnegative integer; hence, $\ell
_{i}\left(  \sigma\right)  \geq0$.
\end{verlong}

We have $i^{\prime}=i+\underbrace{\ell_{i}\left(  \sigma\right)  }_{\geq0}\geq
i$. Hence, the permutation $s_{i^{\prime}-1}\circ s_{i^{\prime}-2}\circ
\cdots\circ s_{i}$ is a composition of $i^{\prime}-i$ simple transpositions
(since $s_{i^{\prime}-1},s_{i^{\prime}-2},\ldots,s_{i}$ are simple
transpositions\footnote{since $s_{k}$ is a simple transposition for each
$k\in\left\{  1,2,\ldots,n-1\right\}  $}). In view of $s_{i^{\prime}-1}\circ
s_{i^{\prime}-2}\circ\cdots\circ s_{i}=a_{i}$ and $\underbrace{i^{\prime}%
}_{=i+\ell_{i}\left(  \sigma\right)  }-i=i+\ell_{i}\left(  \sigma\right)
-i=\ell_{i}\left(  \sigma\right)  $, this rewrites as follows: The permutation
$a_{i}$ is a composition of $\ell_{i}\left(  \sigma\right)  $ simple
transpositions. In other words, $a_{i}=\left(  \text{a composition of }%
\ell_{i}\left(  \sigma\right)  \text{ simple transpositions}\right)  $. This
proves (\ref{sol.perm.lehmer.rothe.d.ai=}).]

Proposition \ref{prop.perm.lehmer.l} yields $\ell\left(  \sigma\right)
=\ell_{1}\left(  \sigma\right)  +\ell_{2}\left(  \sigma\right)  +\cdots
+\ell_{n}\left(  \sigma\right)  $, thus $\ell_{1}\left(  \sigma\right)
+\ell_{2}\left(  \sigma\right)  +\cdots+\ell_{n}\left(  \sigma\right)
=\ell\left(  \sigma\right)  $.

Now, Exercise \ref{exe.perm.lehmer.rothe} \textbf{(c)} yields%
\begin{align*}
\sigma &  =a_{1}\circ a_{2}\circ\cdots\circ a_{n}\\
&  =\left(  \text{a composition of }\ell_{1}\left(  \sigma\right)  \text{
simple transpositions}\right) \\
&  \ \ \ \ \ \ \ \ \ \ \circ\left(  \text{a composition of }\ell_{2}\left(
\sigma\right)  \text{ simple transpositions}\right) \\
&  \ \ \ \ \ \ \ \ \ \ \circ\cdots\\
&  \ \ \ \ \ \ \ \ \ \ \circ\left(  \text{a composition of }\ell_{n}\left(
\sigma\right)  \text{ simple transpositions}\right) \\
&  \ \ \ \ \ \ \ \ \ \ \left(
\begin{array}
[c]{c}%
\text{since }a_{i}=\left(  \text{a composition of }\ell_{i}\left(
\sigma\right)  \text{ simple transpositions}\right) \\
\text{for each }i\in\left\{  1,2,\ldots,n\right\}  \text{ (by
(\ref{sol.perm.lehmer.rothe.d.ai=}))}%
\end{array}
\right) \\
&  =\left(  \text{a composition of }\underbrace{\ell_{1}\left(  \sigma\right)
+\ell_{2}\left(  \sigma\right)  +\cdots+\ell_{n}\left(  \sigma\right)
}_{=\ell\left(  \sigma\right)  }\text{ simple transpositions}\right) \\
&  =\left(  \text{a composition of }\ell\left(  \sigma\right)  \text{ simple
transpositions}\right)  .
\end{align*}
In other words, $\sigma$ is a composition of $\ell\left(  \sigma\right)  $
simple transpositions. In other words, $\sigma$ is a composition of
$\ell\left(  \sigma\right)  $ permutations of the form $s_{k}$ (with
$k\in\left\{  1,2,\ldots,n-1\right\}  $) (because the simple transpositions
are precisely the permutations of the form $s_{k}$ (with $k\in\left\{
1,2,\ldots,n-1\right\}  $)). Thus, Exercise \ref{exe.ps2.2.5} \textbf{(e)} is
solved again.

\textbf{(e)} \textit{Second solution to Exercise \ref{exe.ps2.2.4}
\textbf{(c)}:} Define $\sigma\in S_{n}$ by $\sigma=w_{0}$. For each
$i\in\left[  n\right]  $, we define a permutation $a_{i}\in S_{n}$ as in
Exercise \ref{exe.perm.lehmer.rothe}. Then,%
\begin{equation}
a_{i}=s_{n-1}\circ s_{n-2}\circ\cdots\circ s_{i}\ \ \ \ \ \ \ \ \ \ \text{for
each }i\in\left\{  1,2,\ldots,n\right\}  .
\label{sol.perm.lehmer.rothe.e.liw0.pf.2}%
\end{equation}

[\textit{Proof of (\ref{sol.perm.lehmer.rothe.e.liw0.pf.2}):} Let
$i\in\left\{  1,2,\ldots,n\right\}  $. Thus, $i\in\left\{  1,2,\ldots
,n\right\}  =\left[  n\right]  $.

From $\sigma=w_{0}$, we obtain $\ell_{i}\left(  \sigma\right)  =\ell
_{i}\left(  w_{0}\right)  =n-i$ (by (\ref{sol.perm.lehmer.rothe.e.liw0})). In
other words, $i+\ell_{i}\left(  \sigma\right)  =n$.

Set $i^{\prime}=i+\ell_{i}\left(  \sigma\right)  $. Then, Exercise
\ref{exe.perm.lehmer.rothe} \textbf{(b)} yields
\[
a_{i}=s_{i^{\prime}-1}\circ s_{i^{\prime}-2}\circ\cdots\circ s_{i}%
=s_{n-1}\circ s_{n-2}\circ\cdots\circ s_{i}%
\]
(since $i^{\prime}=i+\ell_{i}\left(  \sigma\right)  =n$). This proves
(\ref{sol.perm.lehmer.rothe.e.liw0.pf.2}).]

Now, from $\sigma=w_{0}$, we get%
\begin{align*}
w_{0}  &  =\sigma=a_{1}\circ a_{2}\circ\cdots\circ a_{n}%
\ \ \ \ \ \ \ \ \ \ \left(  \text{by Exercise \ref{exe.perm.lehmer.rothe}
\textbf{(c)}}\right) \\
&  =\left(  s_{n-1}\circ s_{n-2}\circ\cdots\circ s_{1}\right)  \circ\left(
s_{n-1}\circ s_{n-2}\circ\cdots\circ s_{2}\right)  \circ\cdots\circ\left(
s_{n-1}\circ s_{n-2}\circ\cdots\circ s_{n}\right) \\
&  \ \ \ \ \ \ \ \ \ \ \left(
\begin{array}
[c]{c}%
\text{since }a_{i}=s_{n-1}\circ s_{n-2}\circ\cdots\circ s_{i}\text{ for each
}i\in\left\{  1,2,\ldots,n\right\} \\
\text{(by (\ref{sol.perm.lehmer.rothe.e.liw0.pf.2}))}%
\end{array}
\right)  .
\end{align*}
\footnote{Note that the factor $\left(  s_{n-1}\circ s_{n-2}\circ\cdots\circ
s_{n}\right)  $ is an empty composition of permutations in $S_{n}$, and thus
equals $\operatorname*{id}$. We can omit this factor (when $n>0$).} This is an
explicit way to write $w_{0}$ as a composition of several permutations of the
form $s_{i}$ (with $i\in\left\{  1,2,\ldots,n-1\right\}  $). Thus, Exercise
\ref{exe.ps2.2.4} \textbf{(c)} is solved again.

[\textit{Remark:} The way of writing $w_{0}$ that we have just obtained is
\textbf{not} the exact same way that we found in our first solution to
Exercise \ref{exe.ps2.2.4} \textbf{(c)} above. Indeed, the latter way is%
\[
w_{0}=\left(  s_{0}\circ s_{-1}\circ\cdots\circ s_{1}\right)  \circ\left(
s_{1}\circ s_{0}\circ\cdots\circ s_{1}\right)  \circ\cdots\circ\left(
s_{n-1}\circ s_{n-2}\circ\cdots\circ s_{1}\right)  .
\]
It is, however, possible to derive one of these two representations of $w_{0}$
from the other.]
\end{proof}

\subsection{Solution to Exercise \ref{exe.perm.extend.proofs}}

In this section, we shall use the notations introduced in Definition
\ref{def.perm.extend.SX} and in Definition \ref{def.perm.extend.YtX}.

We begin with the (extremely simple) proof of Proposition
\ref{prop.perm.extend.YtX.wd}:

\begin{vershort}
\begin{proof}
[Proof of Proposition \ref{prop.perm.extend.YtX.wd}.]We just need to prove
that $%
\begin{cases}
\sigma\left(  x\right)  , & \text{if }x\in Y;\\
x, & \text{if }x\notin Y
\end{cases}
\in X$ for each $x\in X$. But this is true, because:

\begin{itemize}
\item If $x\in Y$, then $%
\begin{cases}
\sigma\left(  x\right)  , & \text{if }x\in Y;\\
x, & \text{if }x\notin Y
\end{cases}
=\sigma\left(  x\right)  \in Y\subseteq X$.

\item Otherwise, $%
\begin{cases}
\sigma\left(  x\right)  , & \text{if }x\in Y;\\
x, & \text{if }x\notin Y
\end{cases}
=x\in X$.
\end{itemize}

Thus, Proposition \ref{prop.perm.extend.YtX.wd} is proven.
\end{proof}
\end{vershort}

\begin{verlong}
\begin{proof}
[Proof of Proposition \ref{prop.perm.extend.YtX.wd}.]Every $x\in X$ satisfies
$%
\begin{cases}
\sigma\left(  x\right)  , & \text{if }x\in Y;\\
x, & \text{if }x\notin Y
\end{cases}
\in X$\ \ \ \ \footnote{\textit{Proof.} Fix $x\in X$. We must prove that $%
\begin{cases}
\sigma\left(  x\right)  , & \text{if }x\in Y;\\
x, & \text{if }x\notin Y
\end{cases}
\in X$.
\par
The map $\sigma$ is a map from $Y$ to $Y$; thus, $\sigma\left(  Y\right)
\subseteq Y$.
\par
We are in one of the following two cases:
\par
\textit{Case 1:} We have $x\in Y$.
\par
\textit{Case 2:} We have $x\notin Y$.
\par
Let us first consider Case 1. In this case, we have $x\in Y$. Thus,
\[%
\begin{cases}
\sigma\left(  x\right)  , & \text{if }x\in Y;\\
x, & \text{if }x\notin Y
\end{cases}
=\sigma\left(  \underbrace{x}_{\in Y}\right)  \in\sigma\left(  Y\right)
\subseteq Y\subseteq X.
\]
Hence, $%
\begin{cases}
\sigma\left(  x\right)  , & \text{if }x\in Y;\\
x, & \text{if }x\notin Y
\end{cases}
\in X$ is proven in Case 1.
\par
Let us next consider Case 2. In this case, we have $x\notin Y$. Thus, $%
\begin{cases}
\sigma\left(  x\right)  , & \text{if }x\in Y;\\
x, & \text{if }x\notin Y
\end{cases}
=x\in X$. Hence, $%
\begin{cases}
\sigma\left(  x\right)  , & \text{if }x\in Y;\\
x, & \text{if }x\notin Y
\end{cases}
\in X$ is proven in Case 2.
\par
We have now proven $%
\begin{cases}
\sigma\left(  x\right)  , & \text{if }x\in Y;\\
x, & \text{if }x\notin Y
\end{cases}
\in X$ in each of the two Cases 1 and 2. Since these two Cases cover all
possibilities, we thus conclude that $%
\begin{cases}
\sigma\left(  x\right)  , & \text{if }x\in Y;\\
x, & \text{if }x\notin Y
\end{cases}
\in X$ always holds. Qed.}. Hence, the map $\sigma^{\left(  Y\rightarrow
X\right)  }:X\rightarrow X$ is well-defined\footnote{because it was defined
by
\[
\left(  \sigma^{\left(  Y\rightarrow X\right)  }\left(  x\right)  =%
\begin{cases}
\sigma\left(  x\right)  , & \text{if }x\in Y;\\
x, & \text{if }x\notin Y
\end{cases}
\ \ \ \ \ \ \ \ \ \ \text{for every }x\in X\right)
\]
}. This proves Proposition \ref{prop.perm.extend.YtX.wd}.
\end{proof}
\end{verlong}

Proposition \ref{prop.perm.extend.YtX.wd} justifies Definition
\ref{def.perm.extend.YtX}.

\begin{proof}
[Proof of Proposition \ref{prop.perm.extend.YtX.inj-sur}.]\textbf{(a)} Let
$\alpha:Y\rightarrow Y$ and $\beta:Y\rightarrow Y$ be two maps. Let $x\in X$.
We shall prove that%
\begin{equation}
\left(  \alpha\circ\beta\right)  ^{\left(  Y\rightarrow X\right)  }\left(
x\right)  =\left(  \alpha^{\left(  Y\rightarrow X\right)  }\circ\beta^{\left(
Y\rightarrow X\right)  }\right)  \left(  x\right)  .
\label{pf.prop.perm.extend.YtX.inj-sur.a.goal}%
\end{equation}

[\textit{Proof of (\ref{pf.prop.perm.extend.YtX.inj-sur.a.goal}):} We are in
one of the following two cases:

\textit{Case 1:} We have $x\in Y$.

\textit{Case 2:} We don't have $x\in Y$.

Let us first consider Case 1. In this case, we have $x\in Y$. The definition
of $\beta^{\left(  Y\rightarrow X\right)  }$ yields%
\[
\beta^{\left(  Y\rightarrow X\right)  }\left(  x\right)  =%
\begin{cases}
\beta\left(  x\right)  , & \text{if }x\in Y;\\
x, & \text{if }x\notin Y
\end{cases}
=\beta\left(  x\right)  \ \ \ \ \ \ \ \ \ \ \left(  \text{since }x\in
Y\right)  ,
\]
so that $\beta^{\left(  Y\rightarrow X\right)  }\left(  x\right)
=\beta\left(  \underbrace{x}_{\in Y}\right)  \in\beta\left(  Y\right)
\subseteq Y$ (since $\beta$ is a map from $Y$ to $Y$).

Now,%
\begin{align*}
\left(  \alpha^{\left(  Y\rightarrow X\right)  }\circ\beta^{\left(
Y\rightarrow X\right)  }\right)  \left(  x\right)   &  =\alpha^{\left(
Y\rightarrow X\right)  }\left(  \underbrace{\beta^{\left(  Y\rightarrow
X\right)  }\left(  x\right)  }_{=\beta\left(  x\right)  }\right)
=\alpha^{\left(  Y\rightarrow X\right)  }\left(  \beta\left(  x\right)
\right) \\
&  =%
\begin{cases}
\alpha\left(  \beta\left(  x\right)  \right)  , & \text{if }\beta\left(
x\right)  \in Y;\\
\beta\left(  x\right)  , & \text{if }\beta\left(  x\right)  \notin Y
\end{cases}
\ \ \ \ \ \ \ \ \ \ \left(  \text{by the definition of }\alpha^{\left(
Y\rightarrow X\right)  }\right) \\
&  =\alpha\left(  \beta\left(  x\right)  \right)  \ \ \ \ \ \ \ \ \ \ \left(
\text{since }\beta\left(  x\right)  \in Y\right) \\
&  =\left(  \alpha\circ\beta\right)  \left(  x\right)  .
\end{align*}
Compared with%
\begin{align*}
\left(  \alpha\circ\beta\right)  ^{\left(  Y\rightarrow X\right)  }\left(
x\right)   &  =%
\begin{cases}
\left(  \alpha\circ\beta\right)  \left(  x\right)  , & \text{if }x\in Y;\\
x, & \text{if }x\notin Y
\end{cases}
\ \ \ \ \ \ \ \ \ \ \left(  \text{by the definition of }\left(  \alpha
\circ\beta\right)  ^{\left(  Y\rightarrow X\right)  }\right) \\
&  =\left(  \alpha\circ\beta\right)  \left(  x\right)
\ \ \ \ \ \ \ \ \ \ \left(  \text{since }x\in Y\right)  ,
\end{align*}
this yields $\left(  \alpha\circ\beta\right)  ^{\left(  Y\rightarrow X\right)
}\left(  x\right)  =\left(  \alpha^{\left(  Y\rightarrow X\right)  }\circ
\beta^{\left(  Y\rightarrow X\right)  }\right)  \left(  x\right)  $. We thus
have proven (\ref{pf.prop.perm.extend.YtX.inj-sur.a.goal}) in Case 1.

Let us now consider Case 2. In this case, we don't have $x\in Y$. In other
words, we have $x\notin Y$. The definition of $\beta^{\left(  Y\rightarrow
X\right)  }$ yields%
\[
\beta^{\left(  Y\rightarrow X\right)  }\left(  x\right)  =%
\begin{cases}
\beta\left(  x\right)  , & \text{if }x\in Y;\\
x, & \text{if }x\notin Y
\end{cases}
=x\ \ \ \ \ \ \ \ \ \ \left(  \text{since }x\notin Y\right)  .
\]
Now,%
\begin{align*}
\left(  \alpha^{\left(  Y\rightarrow X\right)  }\circ\beta^{\left(
Y\rightarrow X\right)  }\right)  \left(  x\right)   &  =\alpha^{\left(
Y\rightarrow X\right)  }\left(  \underbrace{\beta^{\left(  Y\rightarrow
X\right)  }\left(  x\right)  }_{=x}\right)  =\alpha^{\left(  Y\rightarrow
X\right)  }\left(  x\right) \\
&  =%
\begin{cases}
\alpha\left(  x\right)  , & \text{if }x\in Y;\\
x, & \text{if }x\notin Y
\end{cases}
\ \ \ \ \ \ \ \ \ \ \left(  \text{by the definition of }\alpha^{\left(
Y\rightarrow X\right)  }\right) \\
&  =x\ \ \ \ \ \ \ \ \ \ \left(  \text{since }x\notin Y\right)  .
\end{align*}
Compared with%
\begin{align*}
\left(  \alpha\circ\beta\right)  ^{\left(  Y\rightarrow X\right)  }\left(
x\right)   &  =%
\begin{cases}
\left(  \alpha\circ\beta\right)  \left(  x\right)  , & \text{if }x\in Y;\\
x, & \text{if }x\notin Y
\end{cases}
\ \ \ \ \ \ \ \ \ \ \left(  \text{by the definition of }\left(  \alpha
\circ\beta\right)  ^{\left(  Y\rightarrow X\right)  }\right) \\
&  =x\ \ \ \ \ \ \ \ \ \ \left(  \text{since }x\notin Y\right)  ,
\end{align*}
this yields $\left(  \alpha\circ\beta\right)  ^{\left(  Y\rightarrow X\right)
}\left(  x\right)  =\left(  \alpha^{\left(  Y\rightarrow X\right)  }\circ
\beta^{\left(  Y\rightarrow X\right)  }\right)  \left(  x\right)  $. We thus
have proven (\ref{pf.prop.perm.extend.YtX.inj-sur.a.goal}) in Case 2.

Now, we have proven the equality (\ref{pf.prop.perm.extend.YtX.inj-sur.a.goal}%
) in each of the two Cases 1 and 2. Since these two Cases cover all
possibilities, this yields that (\ref{pf.prop.perm.extend.YtX.inj-sur.a.goal})
always holds.]

Now, forget that we fixed $x$. We thus have proven that $\left(  \alpha
\circ\beta\right)  ^{\left(  Y\rightarrow X\right)  }\left(  x\right)
=\left(  \alpha^{\left(  Y\rightarrow X\right)  }\circ\beta^{\left(
Y\rightarrow X\right)  }\right)  \left(  x\right)  $ for each $x\in X$. In
other words, $\left(  \alpha\circ\beta\right)  ^{\left(  Y\rightarrow
X\right)  }=\alpha^{\left(  Y\rightarrow X\right)  }\circ\beta^{\left(
Y\rightarrow X\right)  }$. This proves Proposition
\ref{prop.perm.extend.YtX.inj-sur} \textbf{(a)}.

\textbf{(b)} Each $x\in X$ satisfies%
\begin{align*}
\left(  \operatorname*{id}\nolimits_{Y}\right)  ^{\left(  Y\rightarrow
X\right)  }\left(  x\right)   &  =%
\begin{cases}
\operatorname*{id}\nolimits_{Y}\left(  x\right)  , & \text{if }x\in Y;\\
x, & \text{if }x\notin Y
\end{cases}
\ \ \ \ \ \ \ \ \ \ \left(  \text{by the definition of }\left(
\operatorname*{id}\nolimits_{Y}\right)  ^{\left(  Y\rightarrow X\right)
}\right) \\
&  =%
\begin{cases}
x, & \text{if }x\in Y;\\
x, & \text{if }x\notin Y
\end{cases}
\ \ \ \ \ \ \ \ \ \ \left(  \text{since }\operatorname*{id}\nolimits_{Y}%
\left(  x\right)  =x\text{ whenever }x\in Y\right) \\
&  =x=\operatorname*{id}\nolimits_{X}\left(  x\right)  .
\end{align*}
In other words, $\left(  \operatorname*{id}\nolimits_{Y}\right)  ^{\left(
Y\rightarrow X\right)  }=\operatorname*{id}\nolimits_{X}$. This proves
Proposition \ref{prop.perm.extend.YtX.inj-sur} \textbf{(b)}.

\textbf{(c)} Let $\sigma\in S_{Y}$. Thus, $\sigma$ is a permutation of $Y$
(since $S_{Y}$ is the set of all permutations of $Y$). In other words,
$\sigma$ is a bijective map from $Y$ to $Y$. Hence, the map $\sigma^{\left(
Y\rightarrow X\right)  }:X\rightarrow X$ is well-defined. Also, the map
$\sigma$ is bijective, thus invertible. Hence, its inverse $\sigma^{-1}$ is
also a well-defined map from $Y$ to $Y$. Therefore, the map $\left(
\sigma^{-1}\right)  ^{\left(  Y\rightarrow X\right)  }:X\rightarrow X$ is well-defined.

Proposition \ref{prop.perm.extend.YtX.inj-sur} \textbf{(a)} (applied to
$\alpha=\sigma$ and $\beta=\sigma^{-1}$) yields%
\[
\left(  \sigma\circ\sigma^{-1}\right)  ^{\left(  Y\rightarrow X\right)
}=\sigma^{\left(  Y\rightarrow X\right)  }\circ\left(  \sigma^{-1}\right)
^{\left(  Y\rightarrow X\right)  }.
\]
Hence,%
\begin{equation}
\sigma^{\left(  Y\rightarrow X\right)  }\circ\left(  \sigma^{-1}\right)
^{\left(  Y\rightarrow X\right)  }=\left(  \underbrace{\sigma\circ\sigma^{-1}%
}_{=\operatorname*{id}\nolimits_{Y}}\right)  ^{\left(  Y\rightarrow X\right)
}=\left(  \operatorname*{id}\nolimits_{Y}\right)  ^{\left(  Y\rightarrow
X\right)  }=\operatorname*{id}\nolimits_{X}
\label{pf.prop.perm.extend.YtX.inj-sur.c.1}%
\end{equation}
(by Proposition \ref{prop.perm.extend.YtX.inj-sur} \textbf{(b)}). Also,
Proposition \ref{prop.perm.extend.YtX.inj-sur} \textbf{(a)} (applied to
$\alpha=\sigma^{-1}$ and $\beta=\sigma$) yields%
\[
\left(  \sigma^{-1}\circ\sigma\right)  ^{\left(  Y\rightarrow X\right)
}=\left(  \sigma^{-1}\right)  ^{\left(  Y\rightarrow X\right)  }\circ
\sigma^{\left(  Y\rightarrow X\right)  }.
\]
Hence,%
\begin{equation}
\left(  \sigma^{-1}\right)  ^{\left(  Y\rightarrow X\right)  }\circ
\sigma^{\left(  Y\rightarrow X\right)  }=\left(  \underbrace{\sigma^{-1}%
\circ\sigma}_{=\operatorname*{id}\nolimits_{Y}}\right)  ^{\left(  Y\rightarrow
X\right)  }=\left(  \operatorname*{id}\nolimits_{Y}\right)  ^{\left(
Y\rightarrow X\right)  }=\operatorname*{id}\nolimits_{X}
\label{pf.prop.perm.extend.YtX.inj-sur.c.2}%
\end{equation}
(by Proposition \ref{prop.perm.extend.YtX.inj-sur} \textbf{(b)}).

Combining (\ref{pf.prop.perm.extend.YtX.inj-sur.c.1}) with
(\ref{pf.prop.perm.extend.YtX.inj-sur.c.2}), we conclude that the maps
$\sigma^{\left(  Y\rightarrow X\right)  }$ and $\left(  \sigma^{-1}\right)
^{\left(  Y\rightarrow X\right)  }$ are mutually inverse. Thus, the map
$\sigma^{\left(  Y\rightarrow X\right)  }$ is invertible, and hence bijective.
Hence, $\sigma^{\left(  Y\rightarrow X\right)  }$ is a bijective map from $X$
to $X$. In other words, $\sigma^{\left(  Y\rightarrow X\right)  }$ is a
permutation of $X$. In other words, $\sigma^{\left(  Y\rightarrow X\right)
}\in S_{X}$ (since $S_{X}$ is the set of all permutations of $X$).
Furthermore, it satisfies $\left(  \sigma^{-1}\right)  ^{\left(  Y\rightarrow
X\right)  }=\left(  \sigma^{\left(  Y\rightarrow X\right)  }\right)  ^{-1}$
(since the maps $\sigma^{\left(  Y\rightarrow X\right)  }$ and $\left(
\sigma^{-1}\right)  ^{\left(  Y\rightarrow X\right)  }$ are mutually inverse).
Thus, Proposition \ref{prop.perm.extend.YtX.inj-sur} \textbf{(c)} is proven.

\textbf{(d)} First, we claim that%
\begin{align}
&  \left\{  \delta^{\left(  Y\rightarrow X\right)  }\ \mid\ \delta\in
S_{Y}\right\} \nonumber\\
&  \subseteq\left\{  \tau\in S_{X}\ \mid\ \tau\left(  z\right)  =z\text{ for
every }z\in X\setminus Y\right\}  .
\label{pf.prop.perm.extend.YtX.inj-sur.d.1}%
\end{align}

[\textit{Proof of (\ref{pf.prop.perm.extend.YtX.inj-sur.d.1}):} Let $\alpha
\in\left\{  \delta^{\left(  Y\rightarrow X\right)  }\ \mid\ \delta\in
S_{Y}\right\}  $. Thus, $\alpha=\delta^{\left(  Y\rightarrow X\right)  }$ for
some map $\delta\in S_{Y}$. Consider this $\delta$. Proposition
\ref{prop.perm.extend.YtX.inj-sur} \textbf{(c)} (applied to $\sigma=\delta$)
yields that $\delta^{\left(  Y\rightarrow X\right)  }\in S_{X}$ and $\left(
\delta^{-1}\right)  ^{\left(  Y\rightarrow X\right)  }=\left(  \delta^{\left(
Y\rightarrow X\right)  }\right)  ^{-1}$.

Now, for every $z\in X\setminus Y$, we have
\begin{align*}
\underbrace{\alpha}_{=\delta^{\left(  Y\rightarrow X\right)  }}\left(
z\right)   &  =\delta^{\left(  Y\rightarrow X\right)  }\left(  z\right)  =%
\begin{cases}
\delta\left(  z\right)  , & \text{if }z\in Y;\\
z, & \text{if }z\notin Y
\end{cases}
\ \ \ \ \ \ \ \ \ \ \left(  \text{by the definition of }\delta^{\left(
Y\rightarrow X\right)  }\right) \\
&  =z\ \ \ \ \ \ \ \ \ \ \left(  \text{since }z\notin Y\text{ (because }z\in
X\setminus Y\text{)}\right)  .
\end{align*}
Hence, we have shown that $\alpha\left(  z\right)  =z$ for every $z\in
X\setminus Y$. Thus, $\alpha$ is a $\tau\in S_{X}$ that satisfies $\tau\left(
z\right)  =z$ for every $z\in X\setminus Y$ (since $\alpha=\delta^{\left(
Y\rightarrow X\right)  }\in S_{X}$). In other words,%
\[
\alpha\in\left\{  \tau\in S_{X}\ \mid\ \tau\left(  z\right)  =z\text{ for
every }z\in X\setminus Y\right\}  .
\]

\begin{vershort}
Since we have proven this for \textbf{each} $\alpha\in\left\{  \delta^{\left(
Y\rightarrow X\right)  }\ \mid\ \delta\in S_{Y}\right\}  $, we thus conclude
that $\left\{  \delta^{\left(  Y\rightarrow X\right)  }\ \mid\ \delta\in
S_{Y}\right\}  \subseteq\left\{  \tau\in S_{X}\ \mid\ \tau\left(  z\right)
=z\text{ for every }z\in X\setminus Y\right\}  $. This proves
(\ref{pf.prop.perm.extend.YtX.inj-sur.d.1}).]
\end{vershort}

\begin{verlong}
Now, forget that we fixed $\alpha$. Thus, we have proven that \newline%
$\alpha\in\left\{  \tau\in S_{X}\ \mid\ \tau\left(  z\right)  =z\text{ for
every }z\in X\setminus Y\right\}  $ for each\newline$\alpha\in\left\{
\delta^{\left(  Y\rightarrow X\right)  }\ \mid\ \delta\in S_{Y}\right\}  $. In
other words,%
\[
\left\{  \delta^{\left(  Y\rightarrow X\right)  }\ \mid\ \delta\in
S_{Y}\right\}  \subseteq\left\{  \tau\in S_{X}\ \mid\ \tau\left(  z\right)
=z\text{ for every }z\in X\setminus Y\right\}  .
\]
This proves (\ref{pf.prop.perm.extend.YtX.inj-sur.d.1}).]
\end{verlong}

Next, we claim that%
\begin{align}
&  \left\{  \tau\in S_{X}\ \mid\ \tau\left(  z\right)  =z\text{ for every
}z\in X\setminus Y\right\} \nonumber\\
&  \subseteq\left\{  \delta^{\left(  Y\rightarrow X\right)  }\ \mid\ \delta\in
S_{Y}\right\}  . \label{pf.prop.perm.extend.YtX.inj-sur.d.2}%
\end{align}

\begin{vershort}
[\textit{Proof of (\ref{pf.prop.perm.extend.YtX.inj-sur.d.2}):} Let $\eta
\in\left\{  \tau\in S_{X}\ \mid\ \tau\left(  z\right)  =z\text{ for every
}z\in X\setminus Y\right\}  $. Thus, $\eta$ is an element of $S_{X}$ and
satisfies%
\begin{equation}
\eta\left(  z\right)  =z\text{ for every }z\in X\setminus Y.
\label{pf.prop.perm.extend.YtX.inj-sur.d.2.pf.short.fix}%
\end{equation}

We have $\eta\in S_{X}$. In other words, $\eta$ is a permutation of $X$. In
other words, $\eta$ is a bijective map from $X$ to $X$. Hence, the map $\eta$
is bijective, and thus both injective and surjective.

Let $y\in Y$. Assume (for the sake of contradiction) that $\eta\left(
y\right)  \notin Y$. Hence, $\eta\left(  y\right)  \in X\setminus Y$. Thus,
(\ref{pf.prop.perm.extend.YtX.inj-sur.d.2.pf.short.fix}) (applied to
$z=\eta\left(  y\right)  $) yields $\eta\left(  \eta\left(  y\right)  \right)
=\eta\left(  y\right)  $. Since $\eta$ is injective, this yields $\eta\left(
y\right)  =y$, hence $y=\eta\left(  y\right)  \in X\setminus Y$, so that
$y\notin Y$, which contradicts $y\in Y$. Hence, our assumption (that
$\eta\left(  y\right)  \notin Y$) was wrong. Hence, we have $\eta\left(
y\right)  \in Y$.

Now, forget that we fixed $y$. We thus have shown that every $y\in Y$
satisfies $\eta\left(  y\right)  \in Y$. Hence, we can define a map
$\varepsilon:Y\rightarrow Y$ by%
\begin{equation}
\left(  \varepsilon\left(  y\right)  =\eta\left(  y\right)
\ \ \ \ \ \ \ \ \ \ \text{for each }y\in Y\right)  .
\label{pf.prop.perm.extend.YtX.inj-sur.d.2.pf.short.ey=}%
\end{equation}
Consider this $\varepsilon$.

This map $\varepsilon$ is injective\footnote{\textit{Proof.} Let $a\in Y$ and
$b\in Y$ be such that $\varepsilon\left(  a\right)  =\varepsilon\left(
b\right)  $. By the definition of $\varepsilon$, we have $\varepsilon\left(
a\right)  =\eta\left(  a\right)  $ and $\varepsilon\left(  b\right)
=\eta\left(  b\right)  $. Now, $\eta\left(  a\right)  =\varepsilon\left(
a\right)  =\varepsilon\left(  b\right)  =\eta\left(  b\right)  $. Hence, $a=b$
(since the map $\eta$ is injective).
\par
Now, forget that we fixed $a$ and $b$. We thus have shown that every $a\in Y$
and $b\in Y$ satisfying $\varepsilon\left(  a\right)  =\varepsilon\left(
b\right)  $ must satisfy $a=b$. In other words, the map $\varepsilon$ is
injective.} and surjective\footnote{\textit{Proof.} Let $w\in Y$. Then, $w\in
Y\subseteq X=\eta\left(  X\right)  $ (since $\eta$ is surjective). In other
words, there exists some $z\in X$ such that $w=\eta\left(  z\right)  $.
Consider this $z$. If we had $z\in X\setminus Y$, then we would have
\begin{align*}
\eta\left(  z\right)   &  =z\ \ \ \ \ \ \ \ \ \ \left(  \text{by
(\ref{pf.prop.perm.extend.YtX.inj-sur.d.2.pf.short.fix})}\right) \\
&  \notin Y\ \ \ \ \ \ \ \ \ \ \left(  \text{since }z\in X\setminus Y\right)
,
\end{align*}
which would contradict $\eta\left(  z\right)  =w\in Y$. Thus, we cannot have
$z\in X\setminus Y$. Hence, we must have $z\notin X\setminus Y$. Combining
$z\in X$ with $z\notin X\setminus Y$, we obtain $z\in X\setminus\left(
X\setminus Y\right)  \subseteq Y$. Hence, $\varepsilon\left(  z\right)  $ is
well-defined. The definition of $\varepsilon$ yields $\varepsilon\left(
z\right)  =\eta\left(  z\right)  =w$. Hence, $w=\varepsilon\left(
\underbrace{z}_{\in Y}\right)  \in\varepsilon\left(  Y\right)  $.
\par
Now, forget that we fixed $w$. We thus have shown that $w\in\varepsilon\left(
Y\right)  $ for each $w\in Y$. In other words, $Y\subseteq\varepsilon\left(
Y\right)  $. In other words, $\varepsilon$ is surjective.}. Hence, this map
$\varepsilon$ is bijective. Thus, $\varepsilon$ is a bijection from $Y$ to
$Y$. In other words, $\varepsilon$ is a permutation of the set $Y$. In other
words, $\varepsilon\in S_{Y}$.

Now, let $x\in X$. We are going to prove that $\varepsilon^{\left(
Y\rightarrow X\right)  }\left(  x\right)  =\eta\left(  x\right)  $.

Let us rewrite the value $\eta\left(  x\right)  $ depending on whether $x$
belongs to $Y$ or not:

\begin{itemize}
\item If $x\in Y$, then $\eta\left(  x\right)  =\varepsilon\left(  x\right)  $
(because the definition of $\varepsilon$ yields $\varepsilon\left(  x\right)
=\eta\left(  x\right)  $).

\item If $x\notin Y$, then $\eta\left(  x\right)  =x$ (because $x\notin Y$
entails $x\in X\setminus Y$, and therefore
(\ref{pf.prop.perm.extend.YtX.inj-sur.d.2.pf.short.fix}) (applied to $z=x$)
yields $\eta\left(  x\right)  =x$).
\end{itemize}

Combining these two observations, we obtain%
\[
\eta\left(  x\right)  =%
\begin{cases}
\varepsilon\left(  x\right)  , & \text{if }x\in Y;\\
x, & \text{if }x\notin Y
\end{cases}
.
\]
But the definition of $\varepsilon^{\left(  Y\rightarrow X\right)  }$ yields%
\[
\varepsilon^{\left(  Y\rightarrow X\right)  }\left(  x\right)  =%
\begin{cases}
\varepsilon\left(  x\right)  , & \text{if }x\in Y;\\
x, & \text{if }x\notin Y
\end{cases}
.
\]
Comparing the preceding two equations, we obtain $\varepsilon^{\left(
Y\rightarrow X\right)  }\left(  x\right)  =\eta\left(  x\right)  $.

Now, forget that we fixed $x$. We now have proven that $\varepsilon^{\left(
Y\rightarrow X\right)  }\left(  x\right)  =\eta\left(  x\right)  $ for every
$x\in X$. In other words, $\varepsilon^{\left(  Y\rightarrow X\right)  }=\eta
$. Hence, $\eta=\varepsilon^{\left(  Y\rightarrow X\right)  }\in\left\{
\delta^{\left(  Y\rightarrow X\right)  }\ \mid\ \delta\in S_{Y}\right\}  $
(since $\varepsilon\in S_{Y}$).

Now, forget that we fixed $\eta$. We thus have shown that every \newline%
$\eta\in\left\{  \tau\in S_{X}\ \mid\ \tau\left(  z\right)  =z\text{ for every
}z\in X\setminus Y\right\}  $ satisfies $\eta\in\left\{  \delta^{\left(
Y\rightarrow X\right)  }\ \mid\ \delta\in S_{Y}\right\}  $. In other words, we
have%
\begin{equation}
\left\{  \tau\in S_{X}\ \mid\ \tau\left(  z\right)  =z\text{ for every }z\in
X\setminus Y\right\}  \subseteq\left\{  \delta^{\left(  Y\rightarrow X\right)
}\ \mid\ \delta\in S_{Y}\right\}  .\nonumber
\end{equation}
This proves (\ref{pf.prop.perm.extend.YtX.inj-sur.d.2}).]
\end{vershort}

\begin{verlong}
[\textit{Proof of (\ref{pf.prop.perm.extend.YtX.inj-sur.d.2}):} Let $\eta
\in\left\{  \tau\in S_{X}\ \mid\ \tau\left(  z\right)  =z\text{ for every
}z\in X\setminus Y\right\}  $. Thus, $\eta$ is a $\tau\in S_{X}$ satisfying
$\left(  \tau\left(  z\right)  =z\text{ for every }z\in X\setminus Y\right)
$. In other words, $\eta$ is an element of $S_{X}$ and satisfies%
\begin{equation}
\eta\left(  z\right)  =z\text{ for every }z\in X\setminus Y.
\label{pf.prop.perm.extend.YtX.inj-sur.d.2.pf.fix}%
\end{equation}

We have $\eta\in S_{X}$. In other words, $\eta$ is a permutation of $X$ (since
$S_{X}$ is the set of all permutations of $X$). In other words, $\eta$ is a
bijective map from $X$ to $X$. The map $\eta$ is bijective, and thus both
injective and surjective.

Let $y\in Y$. Assume (for the sake of contradiction) that $\eta\left(
y\right)  \in X\setminus Y$. Then,
(\ref{pf.prop.perm.extend.YtX.inj-sur.d.2.pf.fix}) (applied to $z=\eta\left(
y\right)  $) yields $\eta\left(  \eta\left(  y\right)  \right)  =\eta\left(
y\right)  $. Since $\eta$ is injective, this yields $\eta\left(  y\right)
=y$, hence $y=\eta\left(  y\right)  \in X\setminus Y$, so that $y\notin Y$,
which contradicts $y\in Y$. Hence, our assumption (that $\eta\left(  y\right)
\in X\setminus Y$) was wrong. Thus, we don't have $\eta\left(  y\right)  \in
X\setminus Y$. Hence, we have $\eta\left(  y\right)  \notin X\setminus Y$.
Combining $\eta\left(  y\right)  \in X$ with $\eta\left(  y\right)  \notin
X\setminus Y$, we obtain $\eta\left(  y\right)  \in X\setminus\left(
X\setminus Y\right)  \subseteq Y$.

Now, forget that we fixed $y$. We thus have shown that every $y\in Y$
satisfies $\eta\left(  y\right)  \in Y$. Hence, we can define a map
$\varepsilon:Y\rightarrow Y$ by%
\begin{equation}
\left(  \varepsilon\left(  y\right)  =\eta\left(  y\right)
\ \ \ \ \ \ \ \ \ \ \text{for each }y\in Y\right)  .
\label{pf.prop.perm.extend.YtX.inj-sur.d.2.pf.ey=}%
\end{equation}
Consider this $\varepsilon$.

This map $\varepsilon$ is injective\footnote{\textit{Proof.} Let $a\in Y$ and
$b\in Y$ be such that $\varepsilon\left(  a\right)  =\varepsilon\left(
b\right)  $. By the definition of $\varepsilon$, we have $\varepsilon\left(
a\right)  =\eta\left(  a\right)  $ and $\varepsilon\left(  b\right)
=\eta\left(  b\right)  $. Now, $\eta\left(  a\right)  =\varepsilon\left(
a\right)  =\varepsilon\left(  b\right)  =\eta\left(  b\right)  $. Hence, $a=b$
(since the map $\eta$ is injective).
\par
Now, forget that we fixed $a$ and $b$. We thus have shown that every $a\in Y$
and $b\in Y$ satisfying $\varepsilon\left(  a\right)  =\varepsilon\left(
b\right)  $ must satisfy $a=b$. In other words, the map $\varepsilon$ is
injective.} and surjective\footnote{\textit{Proof.} Let $w\in Y$. Then, $w\in
Y\subseteq X=\eta\left(  X\right)  $ (since $\eta$ is surjective). In other
words, there exists some $z\in X$ such that $w=\eta\left(  z\right)  $.
Consider this $z$. If we had $z\in X\setminus Y$, then we would have
\begin{align*}
\eta\left(  z\right)   &  =z\ \ \ \ \ \ \ \ \ \ \left(  \text{by
(\ref{pf.prop.perm.extend.YtX.inj-sur.d.2.pf.fix})}\right) \\
&  \notin Y\ \ \ \ \ \ \ \ \ \ \left(  \text{since }z\in X\setminus Y\right)
,
\end{align*}
which would contradict $\eta\left(  z\right)  =w\in Y$. Thus, we cannot have
$z\in X\setminus Y$. Hence, we must have $z\notin X\setminus Y$. Combining
$z\in X$ with $z\notin X\setminus Y$, we obtain $z\in X\setminus\left(
X\setminus Y\right)  \subseteq Y$. Hence, $\varepsilon\left(  z\right)  $ is
well-defined. The definition of $\varepsilon$ yields $\varepsilon\left(
z\right)  =\eta\left(  z\right)  =w$. Hence, $w=\varepsilon\left(
\underbrace{z}_{\in Y}\right)  \in\varepsilon\left(  Y\right)  $.
\par
Now, forget that we fixed $w$. We thus have shown that $w\in\varepsilon\left(
Y\right)  $ for each $w\in Y$. In other words, $Y\subseteq\varepsilon\left(
Y\right)  $. In other words, $\varepsilon$ is surjective.}. Hence, this map
$\varepsilon$ is bijective. Thus, $\varepsilon$ is a bijection from $Y$ to
$Y$. In other words, $\varepsilon$ is a permutation of the set $Y$. In other
words, $\varepsilon\in S_{Y}$ (since $S_{Y}$ is the set of all permutations of
the set $Y$).

Now, let $x\in X$. We are going to prove that $\varepsilon^{\left(
Y\rightarrow X\right)  }\left(  x\right)  =\eta\left(  x\right)  $.

If $x\in Y$, then%
\begin{align*}
\varepsilon^{\left(  Y\rightarrow X\right)  }\left(  x\right)   &  =%
\begin{cases}
\varepsilon\left(  x\right)  , & \text{if }x\in Y;\\
x, & \text{if }x\notin Y
\end{cases}
\ \ \ \ \ \ \ \ \ \ \left(  \text{by the definition of }\varepsilon^{\left(
Y\rightarrow X\right)  }\right) \\
&  =\varepsilon\left(  x\right)  \ \ \ \ \ \ \ \ \ \ \left(  \text{since }x\in
Y\right) \\
&  =\eta\left(  x\right)  \ \ \ \ \ \ \ \ \ \ \left(  \text{by the definition
of }\varepsilon\right)  .
\end{align*}
Hence, $\varepsilon^{\left(  Y\rightarrow X\right)  }\left(  x\right)
=\eta\left(  x\right)  $ is proven under the assumption that $x\in Y$. Thus,
for the rest of the proof of $\varepsilon^{\left(  Y\rightarrow X\right)
}\left(  x\right)  =\eta\left(  x\right)  $, we can WLOG assume that we don't
have $x\in Y$. Assume this.

We have $x\notin Y$ (since we don't have $x\in Y$). Combining $x\in X$ with
$x\notin Y$, we obtain $x\in X\setminus Y$. Thus, $\eta\left(  x\right)  =x$
(by (\ref{pf.prop.perm.extend.YtX.inj-sur.d.2.pf.fix}), applied to $z=x$).
Now, the definition of $\varepsilon^{\left(  Y\rightarrow X\right)  }$ yields%
\begin{align*}
\varepsilon^{\left(  Y\rightarrow X\right)  }\left(  x\right)   &  =%
\begin{cases}
\varepsilon\left(  x\right)  , & \text{if }x\in Y;\\
x, & \text{if }x\notin Y
\end{cases}
=x\ \ \ \ \ \ \ \ \ \ \left(  \text{since }x\notin Y\right) \\
&  =\eta\left(  x\right)  \ \ \ \ \ \ \ \ \ \ \left(  \text{since }\eta\left(
x\right)  =x\right)  .
\end{align*}
Thus, we have shown that $\varepsilon^{\left(  Y\rightarrow X\right)  }\left(
x\right)  =\eta\left(  x\right)  $.

Now, forget that we fixed $x$. We now have proven that $\varepsilon^{\left(
Y\rightarrow X\right)  }\left(  x\right)  =\eta\left(  x\right)  $ for every
$x\in X$. In other words, $\varepsilon^{\left(  Y\rightarrow X\right)  }=\eta
$. Hence, $\eta=\varepsilon^{\left(  Y\rightarrow X\right)  }$. Thus, the map
$\eta$ has the form $\delta^{\left(  Y\rightarrow X\right)  }$ for some
$\delta\in S_{Y}$ (namely, for $\delta=\varepsilon$), because $\varepsilon\in
S_{Y}$. In other words, $\eta\in\left\{  \delta^{\left(  Y\rightarrow
X\right)  }\ \mid\ \delta\in S_{Y}\right\}  $.

Now, forget that we fixed $\eta$. We thus have shown that every \newline%
$\eta\in\left\{  \tau\in S_{X}\ \mid\ \tau\left(  z\right)  =z\text{ for every
}z\in X\setminus Y\right\}  $ satisfies $\eta\in\left\{  \delta^{\left(
Y\rightarrow X\right)  }\ \mid\ \delta\in S_{Y}\right\}  $. In other words, we
have%
\begin{equation}
\left\{  \tau\in S_{X}\ \mid\ \tau\left(  z\right)  =z\text{ for every }z\in
X\setminus Y\right\}  \subseteq\left\{  \delta^{\left(  Y\rightarrow X\right)
}\ \mid\ \delta\in S_{Y}\right\}  .\nonumber
\end{equation}
This proves (\ref{pf.prop.perm.extend.YtX.inj-sur.d.2}).]
\end{verlong}

Combining the two inclusions (\ref{pf.prop.perm.extend.YtX.inj-sur.d.1}) and
(\ref{pf.prop.perm.extend.YtX.inj-sur.d.2}), we obtain%
\begin{equation}
\left\{  \delta^{\left(  Y\rightarrow X\right)  }\ \mid\ \delta\in
S_{Y}\right\}  =\left\{  \tau\in S_{X}\ \mid\ \tau\left(  z\right)  =z\text{
for every }z\in X\setminus Y\right\}  .\nonumber
\end{equation}
This proves Proposition \ref{prop.perm.extend.YtX.inj-sur} \textbf{(d)}.

\textbf{(e)} Define a subset $G$ of $S_{X}$ by
\begin{equation}
G=\left\{  \tau\in S_{X}\ \mid\ \tau\left(  z\right)  =z\text{ for every }z\in
X\setminus Y\right\}  . \label{pf.prop.perm.extend.YtX.inj-sur.e.G=}%
\end{equation}

Proposition \ref{prop.perm.extend.YtX.inj-sur} \textbf{(d)} yields%
\begin{align}
\left\{  \delta^{\left(  Y\rightarrow X\right)  }\ \mid\ \delta\in
S_{Y}\right\}   &  =\left\{  \tau\in S_{X}\ \mid\ \tau\left(  z\right)
=z\text{ for every }z\in X\setminus Y\right\} \nonumber\\
&  =G\ \ \ \ \ \ \ \ \ \ \left(  \text{by
(\ref{pf.prop.perm.extend.YtX.inj-sur.e.G=})}\right)  .
\label{pf.prop.perm.extend.YtX.inj-sur.e.del=G}%
\end{align}
Thus, $\left\{  \delta^{\left(  Y\rightarrow X\right)  }\ \mid\ \delta\in
S_{Y}\right\}  =G\subseteq G$. In other words, we have $\delta^{\left(
Y\rightarrow X\right)  }\in G$ for each $\delta\in S_{Y}$. Hence, we can
define a map $R:S_{Y}\rightarrow G$ by%
\[
\left(  R\left(  \delta\right)  =\delta^{\left(  Y\rightarrow X\right)
}\ \ \ \ \ \ \ \ \ \ \text{for each }\delta\in S_{Y}\right)  .
\]
Consider this map $R$.

The map $R$ is injective\footnote{\textit{Proof.} Let $\delta\in S_{Y}$ and
$\varepsilon\in S_{Y}$ be such that $R\left(  \delta\right)  =R\left(
\varepsilon\right)  $. We shall prove that $\delta=\varepsilon$.
\par
Let $y\in Y$. Thus, $y\in Y\subseteq X$. The definition of $R$ yields
$R\left(  \delta\right)  =\delta^{\left(  Y\rightarrow X\right)  }$ and
$R\left(  \varepsilon\right)  =\varepsilon^{\left(  Y\rightarrow X\right)  }$.
Hence, $\delta^{\left(  Y\rightarrow X\right)  }=R\left(  \delta\right)
=R\left(  \varepsilon\right)  =\varepsilon^{\left(  Y\rightarrow X\right)  }$.
\par
But $y\in X$. Thus, the definition of $\delta^{\left(  Y\rightarrow X\right)
}$ yields%
\[
\delta^{\left(  Y\rightarrow X\right)  }\left(  y\right)  =%
\begin{cases}
\delta\left(  y\right)  , & \text{if }y\in Y;\\
y, & \text{if }y\notin Y
\end{cases}
=\delta\left(  y\right)  \ \ \ \ \ \ \ \ \ \ \left(  \text{since }y\in
Y\right)  .
\]
The same argument (applied to $\varepsilon$ instead of $\delta$) yields
$\varepsilon^{\left(  Y\rightarrow X\right)  }\left(  y\right)  =\varepsilon
\left(  y\right)  $. Now, from $\delta^{\left(  Y\rightarrow X\right)
}=\varepsilon^{\left(  Y\rightarrow X\right)  }$, we obtain $\delta^{\left(
Y\rightarrow X\right)  }\left(  y\right)  =\varepsilon^{\left(  Y\rightarrow
X\right)  }\left(  y\right)  =\varepsilon\left(  y\right)  $. Comparing this
with $\delta^{\left(  Y\rightarrow X\right)  }\left(  y\right)  =\delta\left(
y\right)  $, we obtain $\delta\left(  y\right)  =\varepsilon\left(  y\right)
$.
\par
Now, forget that we fixed $y$. We thus have shown that $\delta\left(
y\right)  =\varepsilon\left(  y\right)  $ for each $y\in Y$. In other words,
$\delta=\varepsilon$.
\par
Now, forget that we fixed $\delta$ and $\varepsilon$. We thus have proven that
if $\delta\in S_{Y}$ and $\varepsilon\in S_{Y}$ are such that $R\left(
\delta\right)  =R\left(  \varepsilon\right)  $, then $\delta=\varepsilon$. In
other words, the map $R$ is injective.} and
surjective\footnote{\textit{Proof.} Let $\eta\in G$. Then,
\[
\eta\in G=\left\{  \delta^{\left(  Y\rightarrow X\right)  }\ \mid\ \delta\in
S_{Y}\right\}  \ \ \ \ \ \ \ \ \ \ \left(  \text{by
(\ref{pf.prop.perm.extend.YtX.inj-sur.e.G=})}\right)  .
\]
In other words, there exists some $\delta\in S_{Y}$ such that $\eta
=\delta^{\left(  Y\rightarrow X\right)  }$. Consider this $\delta$.
\par
Now, the definition of $R\left(  \delta\right)  $ yields $R\left(
\delta\right)  =\delta^{\left(  Y\rightarrow X\right)  }$. Comparing this with
$\eta=\delta^{\left(  Y\rightarrow X\right)  }$, we obtain $\eta=R\left(
\underbrace{\delta}_{\in S_{Y}}\right)  \in R\left(  S_{Y}\right)  $.
\par
Now, forget that we fixed $\eta$. We thus have shown that $\eta\in R\left(
S_{Y}\right)  $ for each $\eta\in G$. In other words, $G\subseteq R\left(
S_{Y}\right)  $. In other words, the map $R$ is surjective.}. Hence, the map
$R$ is bijective.

Now, the map $R$ is a map from $S_{Y}$ to $G$ such that $\left(  R\left(
\delta\right)  =\delta^{\left(  Y\rightarrow X\right)  }\text{ for each
}\delta\in S_{Y}\right)  $. In view of
(\ref{pf.prop.perm.extend.YtX.inj-sur.e.G=}), this rewrites as follows: The
map $R$ is a map from $S_{Y}$ to \newline$\left\{  \tau\in S_{X}\ \mid
\ \tau\left(  z\right)  =z\text{ for every }z\in X\setminus Y\right\}  $ such
that $\left(  R\left(  \delta\right)  =\delta^{\left(  Y\rightarrow X\right)
}\text{ for each }\delta\in S_{Y}\right)  $. Hence, $R$ is the map
\begin{align*}
S_{Y}  &  \rightarrow\left\{  \tau\in S_{X}\ \mid\ \tau\left(  z\right)
=z\text{ for every }z\in X\setminus Y\right\}  ,\\
\delta &  \mapsto\delta^{\left(  Y\rightarrow X\right)  }.
\end{align*}
Thus, the map%
\begin{align*}
S_{Y}  &  \rightarrow\left\{  \tau\in S_{X}\ \mid\ \tau\left(  z\right)
=z\text{ for every }z\in X\setminus Y\right\}  ,\\
\delta &  \mapsto\delta^{\left(  Y\rightarrow X\right)  }%
\end{align*}
is well-defined and bijective (since $R$ is bijective). This proves
Proposition \ref{prop.perm.extend.YtX.inj-sur} \textbf{(e)}.
\end{proof}

\begin{proof}
[Proof of Proposition \ref{prop.perm.extend.YtX.ZtY}.]Let $x\in X$. We are
going to show that $\left(  \sigma^{\left(  Z\rightarrow Y\right)  }\right)
^{\left(  Y\rightarrow X\right)  }\left(  x\right)  =\sigma^{\left(
Z\rightarrow X\right)  }\left(  x\right)  $.

The definition of $\sigma^{\left(  Z\rightarrow X\right)  }$ yields%
\begin{equation}
\sigma^{\left(  Z\rightarrow X\right)  }\left(  x\right)  =%
\begin{cases}
\sigma\left(  x\right)  , & \text{if }x\in Z;\\
x, & \text{if }x\notin Z
\end{cases}
. \label{pf.prop.perm.extend.YtX.ZtY.6}%
\end{equation}

Now, we distinguish between two cases:

\textit{Case 1:} We have $x\in Y$.

\textit{Case 2:} We don't have $x\in Y$.

Let us first consider Case 1. In this case, we have $x\in Y$. The definition
of $\left(  \sigma^{\left(  Z\rightarrow Y\right)  }\right)  ^{\left(
Y\rightarrow X\right)  }$ yields
\begin{align*}
\left(  \sigma^{\left(  Z\rightarrow Y\right)  }\right)  ^{\left(
Y\rightarrow X\right)  }\left(  x\right)   &  =%
\begin{cases}
\sigma^{\left(  Z\rightarrow Y\right)  }\left(  x\right)  , & \text{if }x\in
Y;\\
x, & \text{if }x\notin Y
\end{cases}
=\sigma^{\left(  Z\rightarrow Y\right)  }\left(  x\right)
\ \ \ \ \ \ \ \ \ \ \left(  \text{since }x\in Y\right) \\
&  =%
\begin{cases}
\sigma\left(  x\right)  , & \text{if }x\in Z;\\
x, & \text{if }x\notin Z
\end{cases}
\ \ \ \ \ \ \ \ \ \ \left(  \text{by the definition of }\sigma^{\left(
Z\rightarrow Y\right)  }\right) \\
&  =\sigma^{\left(  Z\rightarrow X\right)  }\left(  x\right)
\ \ \ \ \ \ \ \ \ \ \left(  \text{by (\ref{pf.prop.perm.extend.YtX.ZtY.6}%
)}\right)  .
\end{align*}
We thus have shown that $\left(  \sigma^{\left(  Z\rightarrow Y\right)
}\right)  ^{\left(  Y\rightarrow X\right)  }\left(  x\right)  =\sigma^{\left(
Z\rightarrow X\right)  }\left(  x\right)  $ in Case 1.

Let us now consider Case 2. In this case, we don't have $x\in Y$. In other
words, we have $x\notin Y$. Thus, $x\notin Z$ (since otherwise, we would have
$x\in Z\subseteq Y$, contradicting $x\notin Y$). Thus,
(\ref{pf.prop.perm.extend.YtX.ZtY.6}) becomes%
\begin{equation}
\sigma^{\left(  Z\rightarrow X\right)  }\left(  x\right)  =%
\begin{cases}
\sigma\left(  x\right)  , & \text{if }x\in Z;\\
x, & \text{if }x\notin Z
\end{cases}
=x\ \ \ \ \ \ \ \ \ \ \left(  \text{since }x\notin Z\right)  .
\label{pf.prop.perm.extend.YtX.ZtY.8}%
\end{equation}

Now, the definition of $\left(  \sigma^{\left(  Z\rightarrow Y\right)
}\right)  ^{\left(  Y\rightarrow X\right)  }$ yields
\begin{align*}
\left(  \sigma^{\left(  Z\rightarrow Y\right)  }\right)  ^{\left(
Y\rightarrow X\right)  }\left(  x\right)   &  =%
\begin{cases}
\sigma^{\left(  Z\rightarrow Y\right)  }\left(  x\right)  , & \text{if }x\in
Y;\\
x, & \text{if }x\notin Y
\end{cases}
=x\ \ \ \ \ \ \ \ \ \ \left(  \text{since }x\notin Y\right) \\
&  =\sigma^{\left(  Z\rightarrow X\right)  }\left(  x\right)
\ \ \ \ \ \ \ \ \ \ \left(  \text{by (\ref{pf.prop.perm.extend.YtX.ZtY.8}%
)}\right)  .
\end{align*}
We thus have shown that $\left(  \sigma^{\left(  Z\rightarrow Y\right)
}\right)  ^{\left(  Y\rightarrow X\right)  }\left(  x\right)  =\sigma^{\left(
Z\rightarrow X\right)  }\left(  x\right)  $ in Case 2.

We now have proven that $\left(  \sigma^{\left(  Z\rightarrow Y\right)
}\right)  ^{\left(  Y\rightarrow X\right)  }\left(  x\right)  =\sigma^{\left(
Z\rightarrow X\right)  }\left(  x\right)  $ holds in each of the two Cases 1
and 2. Since these two Cases cover all possibilities, this yields that
$\left(  \sigma^{\left(  Z\rightarrow Y\right)  }\right)  ^{\left(
Y\rightarrow X\right)  }\left(  x\right)  =\sigma^{\left(  Z\rightarrow
X\right)  }\left(  x\right)  $ always holds.

Now, forget that we fixed $x$. We thus have shown that $\left(  \sigma
^{\left(  Z\rightarrow Y\right)  }\right)  ^{\left(  Y\rightarrow X\right)
}\left(  x\right)  =\sigma^{\left(  Z\rightarrow X\right)  }\left(  x\right)
$ for every $x\in X$. In other words, $\left(  \sigma^{\left(  Z\rightarrow
Y\right)  }\right)  ^{\left(  Y\rightarrow X\right)  }=\sigma^{\left(
Z\rightarrow X\right)  }$. This proves Proposition
\ref{prop.perm.extend.YtX.ZtY}.
\end{proof}

\begin{proof}
[Proof of Proposition \ref{prop.perm.extend.YtX.commute}.]We have
$X\setminus\left(  X\setminus Y\right)  =Y$ (since $Y\subseteq X$).

Define a map $\eta:X\rightarrow X$ by $\eta=\alpha^{\left(  Y\rightarrow
X\right)  }$.

Define a map $\zeta:X\rightarrow X$ by $\zeta=\beta^{\left(  X\setminus
Y\rightarrow X\right)  }$.

Let now $x\in X$. We are going to prove that $\left(  \eta\circ\zeta\right)
\left(  x\right)  =\left(  \zeta\circ\eta\right)  \left(  x\right)  $.

We distinguish between two cases:

\textit{Case 1:} We have $x\in Y$.

\textit{Case 2:} We don't have $x\in Y$.

Let us first consider Case 1. In this case, we have $x\in Y$. Thus, $x\in
Y=X\setminus\left(  X\setminus Y\right)  $. Hence, $x\in X$ and $x\notin
X\setminus Y$. Now, recall that $\zeta=\beta^{\left(  X\setminus Y\rightarrow
X\right)  }$. Hence,%
\begin{align*}
\zeta\left(  x\right)   &  =\beta^{\left(  X\setminus Y\rightarrow X\right)
}\left(  x\right)  =%
\begin{cases}
\beta\left(  x\right)  , & \text{if }x\in X\setminus Y;\\
x, & \text{if }x\notin X\setminus Y
\end{cases}
\ \ \ \ \ \ \ \ \ \ \left(  \text{by the definition of }\beta^{\left(
X\setminus Y\rightarrow X\right)  }\right) \\
&  =x\ \ \ \ \ \ \ \ \ \ \left(  \text{since }x\notin X\setminus Y\right)  .
\end{align*}
Also, $\eta=\alpha^{\left(  Y\rightarrow X\right)  }$, so that%
\begin{align*}
\eta\left(  x\right)   &  =\alpha^{\left(  Y\rightarrow X\right)  }\left(
x\right)  =%
\begin{cases}
\alpha\left(  x\right)  , & \text{if }x\in Y;\\
x, & \text{if }x\notin Y
\end{cases}
\ \ \ \ \ \ \ \ \ \ \left(  \text{by the definition of }\alpha^{\left(
Y\rightarrow X\right)  }\right) \\
&  =\alpha\left(  \underbrace{x}_{\in Y}\right)  \ \ \ \ \ \ \ \ \ \ \left(
\text{since }x\in Y\right) \\
&  \in\alpha\left(  Y\right)  \subseteq Y\ \ \ \ \ \ \ \ \ \ \left(
\text{since }\alpha\text{ is a map }Y\rightarrow Y\right) \\
&  =X\setminus\left(  X\setminus Y\right)  ,
\end{align*}
and thus $\eta\left(  x\right)  \notin X\setminus Y$. Now,%
\[
\left(  \eta\circ\zeta\right)  \left(  x\right)  =\eta\left(
\underbrace{\zeta\left(  x\right)  }_{=x}\right)  =\eta\left(  x\right)  .
\]
Compared with%
\begin{align*}
\left(  \zeta\circ\eta\right)  \left(  x\right)   &  =\underbrace{\zeta
}_{=\beta^{\left(  X\setminus Y\rightarrow X\right)  }}\left(  \eta\left(
x\right)  \right)  =\beta^{\left(  X\setminus Y\rightarrow X\right)  }\left(
\eta\left(  x\right)  \right) \\
&  =%
\begin{cases}
\beta\left(  \eta\left(  x\right)  \right)  , & \text{if }\eta\left(
x\right)  \in X\setminus Y;\\
\eta\left(  x\right)  , & \text{if }\eta\left(  x\right)  \notin X\setminus Y
\end{cases}
\ \ \ \ \ \ \ \ \ \ \left(  \text{by the definition of }\beta^{\left(
X\setminus Y\rightarrow X\right)  }\right) \\
&  =\eta\left(  x\right)  \ \ \ \ \ \ \ \ \ \ \left(  \text{since }\eta\left(
x\right)  \notin X\setminus Y\right)  ,
\end{align*}
this yields $\left(  \eta\circ\zeta\right)  \left(  x\right)  =\left(
\zeta\circ\eta\right)  \left(  x\right)  $. Thus, $\left(  \eta\circ
\zeta\right)  \left(  x\right)  =\left(  \zeta\circ\eta\right)  \left(
x\right)  $ is proven in Case 1.

Now, let us consider Case 2. In this case, we don't have $x\in Y$. Thus,
$x\notin Y$. Combined with $x\in X$, this yields $x\in X\setminus Y$. Now,
recall that $\zeta=\beta^{\left(  X\setminus Y\rightarrow X\right)  }$, so
that
\begin{align*}
\zeta\left(  x\right)   &  =\beta^{\left(  X\setminus Y\rightarrow X\right)
}\left(  x\right) \\
&  =%
\begin{cases}
\beta\left(  x\right)  , & \text{if }x\in X\setminus Y;\\
x, & \text{if }x\notin X\setminus Y
\end{cases}
\ \ \ \ \ \ \ \ \ \ \left(  \text{by the definition of }\beta^{\left(
X\setminus Y\rightarrow X\right)  }\right) \\
&  =\beta\left(  \underbrace{x}_{\in X\setminus Y}\right)
\ \ \ \ \ \ \ \ \ \ \left(  \text{since }x\in X\setminus Y\right) \\
&  \in\beta\left(  X\setminus Y\right)  \subseteq X\setminus
Y\ \ \ \ \ \ \ \ \ \ \left(  \text{since }\beta\text{ is a map }X\setminus
Y\rightarrow X\setminus Y\right)  ,
\end{align*}
and thus $\zeta\left(  x\right)  \notin Y$. Also, $\eta=\alpha^{\left(
Y\rightarrow X\right)  }$, so that%
\begin{align*}
\eta\left(  x\right)   &  =\alpha^{\left(  Y\rightarrow X\right)  }\left(
x\right)  =%
\begin{cases}
\alpha\left(  x\right)  , & \text{if }x\in Y;\\
x, & \text{if }x\notin Y
\end{cases}
\ \ \ \ \ \ \ \ \ \ \left(  \text{by the definition of }\alpha^{\left(
Y\rightarrow X\right)  }\right) \\
&  =x\ \ \ \ \ \ \ \ \ \ \left(  \text{since }x\notin Y\right)  .
\end{align*}
Now,%
\begin{align*}
\left(  \eta\circ\zeta\right)  \left(  x\right)   &  =\eta\left(  \zeta\left(
x\right)  \right)  =\alpha^{\left(  Y\rightarrow X\right)  }\left(
\zeta\left(  x\right)  \right)  \ \ \ \ \ \ \ \ \ \ \left(  \text{since }%
\eta=\alpha^{\left(  Y\rightarrow X\right)  }\right) \\
&  =%
\begin{cases}
\alpha\left(  \zeta\left(  x\right)  \right)  , & \text{if }\zeta\left(
x\right)  \in Y;\\
\zeta\left(  x\right)  , & \text{if }\zeta\left(  x\right)  \notin Y
\end{cases}
\ \ \ \ \ \ \ \ \ \ \left(  \text{by the definition of }\alpha^{\left(
Y\rightarrow X\right)  }\right) \\
&  =\zeta\left(  x\right)  \ \ \ \ \ \ \ \ \ \ \left(  \text{since }%
\zeta\left(  x\right)  \notin Y\right)  .
\end{align*}
Compared with%
\[
\left(  \zeta\circ\eta\right)  \left(  x\right)  =\zeta\left(
\underbrace{\eta\left(  x\right)  }_{=x}\right)  =\zeta\left(  x\right)  ,
\]
this yields $\left(  \eta\circ\zeta\right)  \left(  x\right)  =\left(
\zeta\circ\eta\right)  \left(  x\right)  $. Thus, $\left(  \eta\circ
\zeta\right)  \left(  x\right)  =\left(  \zeta\circ\eta\right)  \left(
x\right)  $ is proven in Case 2.

We now have proven the equality $\left(  \eta\circ\zeta\right)  \left(
x\right)  =\left(  \zeta\circ\eta\right)  \left(  x\right)  $ in each of the
two Cases 1 and 2. Since these two Cases cover all possibilities, this yields
that the equality $\left(  \eta\circ\zeta\right)  \left(  x\right)  =\left(
\zeta\circ\eta\right)  \left(  x\right)  $ always holds.

Now, forget that we fixed $x$. We thus have proven that $\left(  \eta
\circ\zeta\right)  \left(  x\right)  =\left(  \zeta\circ\eta\right)  \left(
x\right)  $ for every $x\in X$. In other words, $\eta\circ\zeta=\zeta\circ
\eta$. In view of $\eta=\alpha^{\left(  Y\rightarrow X\right)  }$ and
$\zeta=\beta^{\left(  X\setminus Y\rightarrow X\right)  }$, this rewrites as
\[
\alpha^{\left(  Y\rightarrow X\right)  }\circ\beta^{\left(  X\setminus
Y\rightarrow X\right)  }=\beta^{\left(  X\setminus Y\rightarrow X\right)
}\circ\alpha^{\left(  Y\rightarrow X\right)  }.
\]
This proves Proposition \ref{prop.perm.extend.YtX.commute}.
\end{proof}

\begin{proof}
[Proof of Proposition \ref{prop.perm.extend.conj-inv.1}.]We will prove
Proposition \ref{prop.perm.extend.conj-inv.1} by strong induction over
$\left\vert X\right\vert $:

\textit{Induction step:} Let $N\in\mathbb{N}$. Assume that Proposition
\ref{prop.perm.extend.conj-inv.1} holds whenever $\left\vert X\right\vert <N$.
We need to prove that Proposition \ref{prop.perm.extend.conj-inv.1} holds
whenever $\left\vert X\right\vert =N$.

We have assumed that Proposition \ref{prop.perm.extend.conj-inv.1} holds
whenever $\left\vert X\right\vert <N$. In other words, the following claim holds:

\begin{statement}
\textit{Claim 1:} Let $X$ be a finite set satisfying $\left\vert X\right\vert
<N$. Let $x\in X$. Let $\pi\in S_{X}$. Then, there exists a $\sigma\in\left\{
\delta^{\left(  X\setminus\left\{  x\right\}  \rightarrow X\right)  }%
\ \mid\ \delta\in S_{X\setminus\left\{  x\right\}  }\right\}  $ such that
$\sigma\circ\pi\circ\sigma^{-1}=\pi^{-1}$.
\end{statement}

We need to prove that Proposition \ref{prop.perm.extend.conj-inv.1} holds
whenever $\left\vert X\right\vert =N$. In other words, we need to prove the
following claim:

\begin{statement}
\textit{Claim 2:} Let $X$ be a finite set satisfying $\left\vert X\right\vert
=N$. Let $x\in X$. Let $\pi\in S_{X}$. Then, there exists a $\sigma\in\left\{
\delta^{\left(  X\setminus\left\{  x\right\}  \rightarrow X\right)  }%
\ \mid\ \delta\in S_{X\setminus\left\{  x\right\}  }\right\}  $ such that
$\sigma\circ\pi\circ\sigma^{-1}=\pi^{-1}$.
\end{statement}

Before we start proving Claim 2, let us first establish an auxiliary claim,
which follows easily from Claim 1:

\begin{statement}
\textit{Claim 3:} Let $X$ be a finite set satisfying $\left\vert X\right\vert
<N$. Let $\pi\in S_{X}$. Then, there exists a $\sigma\in S_{X}$ such that
$\sigma\circ\pi\circ\sigma^{-1}=\pi^{-1}$.
\end{statement}

[\textit{Proof of Claim 3:} If $X=\varnothing$, then Claim 3 holds for simple
reasons\footnote{\textit{Proof.} Assume that $X=\varnothing$. We must prove
that Claim 3 holds.
\par
We have $X=\varnothing$. Thus, no $x\in X$ exists. Hence, $\left(  \pi\circ
\pi\circ\pi^{-1}\right)  \left(  x\right)  =\pi^{-1}\left(  x\right)  $ holds
for all $x\in X$ (indeed, this is vacuously true, since no $x\in X$ exists).
In other words, we have $\pi\circ\pi\circ\pi^{-1}=\pi^{-1}$. Thus, there
exists a $\sigma\in S_{X}$ such that $\sigma\circ\pi\circ\sigma^{-1}=\pi^{-1}$
(namely, $\sigma=\pi$). Hence, Claim 3 holds.
\par
We thus have proven that if $X=\varnothing$, then Claim 3 holds. Qed.}. Hence,
for the rest of this proof of Claim 3, we can WLOG assume that we don't have
$X=\varnothing$. Assume this.

\begin{vershort}
There exists some $x\in X$ (since we don't have $X=\varnothing$). Consider
this $x$.
\end{vershort}

\begin{verlong}
We don't have $X=\varnothing$. Thus, the set $X$ is nonempty. Hence, there
exists some $x\in X$. Consider this $x$.
\end{verlong}

Claim 1 shows that there exists a $\sigma\in\left\{  \delta^{\left(
X\setminus\left\{  x\right\}  \rightarrow X\right)  }\ \mid\ \delta\in
S_{X\setminus\left\{  x\right\}  }\right\}  $ such that $\sigma\circ\pi
\circ\sigma^{-1}=\pi^{-1}$. Consider this $\sigma$, and denote it by $\alpha$.
Thus, $\alpha$ is an element of $\left\{  \delta^{\left(  X\setminus\left\{
x\right\}  \rightarrow X\right)  }\ \mid\ \delta\in S_{X\setminus\left\{
x\right\}  }\right\}  $ and satisfies $\alpha\circ\pi\circ\alpha^{-1}=\pi
^{-1}$.

We know that $\alpha$ is an element of $\left\{  \delta^{\left(
X\setminus\left\{  x\right\}  \rightarrow X\right)  }\ \mid\ \delta\in
S_{X\setminus\left\{  x\right\}  }\right\}  $. In other words, there exists
some $\delta\in S_{X\setminus\left\{  x\right\}  }$ such that $\alpha
=\delta^{\left(  X\setminus\left\{  x\right\}  \rightarrow X\right)  }$.
Consider this $\delta$.

Now, Proposition \ref{prop.perm.extend.YtX.inj-sur} \textbf{(c)} (applied to
$Y=X\setminus\left\{  x\right\}  $ and $\sigma=\delta$) yields that $\delta$
satisfies $\delta^{\left(  X\setminus\left\{  x\right\}  \rightarrow X\right)
}\in S_{X}$ and $\left(  \delta^{-1}\right)  ^{\left(  X\setminus\left\{
x\right\}  \rightarrow X\right)  }=\left(  \delta^{\left(  X\setminus\left\{
x\right\}  \rightarrow X\right)  }\right)  ^{-1}$. Now, $\alpha=\delta
^{\left(  X\setminus\left\{  x\right\}  \rightarrow X\right)  }\in S_{X}$ and
$\alpha\circ\pi\circ\alpha^{-1}=\pi^{-1}$. Hence, there exists a $\sigma\in
S_{X}$ such that $\sigma\circ\pi\circ\sigma^{-1}=\pi^{-1}$ (namely,
$\sigma=\alpha$). This proves Claim 3.]

For the sake of convenience, we shall also derive a simple consequence of
Proposition \ref{prop.perm.extend.YtX.inj-sur}:

\begin{statement}
\textit{Claim 4:} Let $X$ be a set. Let $x\in X$. Let $\gamma\in
S_{X\setminus\left\{  x\right\}  }$ and $\varepsilon\in S_{X\setminus\left\{
x\right\}  }$. Then:

\textbf{(a)} We have $\gamma^{\left(  X\setminus\left\{  x\right\}
\rightarrow X\right)  }\in S_{X}$.

\textbf{(b)} We have $\left(  \gamma\circ\varepsilon\right)  ^{\left(
X\setminus\left\{  x\right\}  \rightarrow X\right)  }=\gamma^{\left(
X\setminus\left\{  x\right\}  \rightarrow X\right)  }\circ\varepsilon^{\left(
X\setminus\left\{  x\right\}  \rightarrow X\right)  }$.

\textbf{(c)} We have%
\[
\left(  \gamma\circ\varepsilon\circ\gamma^{-1}\right)  ^{\left(
X\setminus\left\{  x\right\}  \rightarrow X\right)  }=\gamma^{\left(
X\setminus\left\{  x\right\}  \rightarrow X\right)  }\circ\varepsilon^{\left(
X\setminus\left\{  x\right\}  \rightarrow X\right)  }\circ\left(
\gamma^{\left(  X\setminus\left\{  x\right\}  \rightarrow X\right)  }\right)
^{-1}.
\]

\end{statement}

[\textit{Proof of Claim 4:} Clearly, $\gamma$ and $\varepsilon$ are elements
of $S_{X\setminus\left\{  x\right\}  }$, thus permutations of the set
$X\setminus\left\{  x\right\}  $ (since $S_{X\setminus\left\{  x\right\}  }$
is the set of all permutations of the set $X\setminus\left\{  x\right\}  $).
In other words, $\gamma$ and $\varepsilon$ are bijective maps from
$X\setminus\left\{  x\right\}  $ to $X\setminus\left\{  x\right\}  $. Thus,
their inverses $\gamma^{-1}$ and $\varepsilon^{-1}$ are maps from
$X\setminus\left\{  x\right\}  $ to $X\setminus\left\{  x\right\}  $ as well.

Proposition \ref{prop.perm.extend.YtX.inj-sur} \textbf{(c)} (applied to
$Y=X\setminus\left\{  x\right\}  $ and $\sigma=\gamma$) shows that $\gamma$
satisfies $\gamma^{\left(  X\setminus\left\{  x\right\}  \rightarrow X\right)
}\in S_{X}$ and $\left(  \gamma^{-1}\right)  ^{\left(  X\setminus\left\{
x\right\}  \rightarrow X\right)  }=\left(  \gamma^{\left(  X\setminus\left\{
x\right\}  \rightarrow X\right)  }\right)  ^{-1}$. This proves Claim 4
\textbf{(a)}. Also, Proposition \ref{prop.perm.extend.YtX.inj-sur}
\textbf{(a)} (applied to $X\setminus\left\{  x\right\}  $, $\gamma$ and
$\varepsilon$ instead of $Y$, $\alpha$ and $\beta$) yields $\left(
\gamma\circ\varepsilon\right)  ^{\left(  X\setminus\left\{  x\right\}
\rightarrow X\right)  }=\gamma^{\left(  X\setminus\left\{  x\right\}
\rightarrow X\right)  }\circ\varepsilon^{\left(  X\setminus\left\{  x\right\}
\rightarrow X\right)  }$. This proves Claim 4 \textbf{(b)}. Now, Proposition
\ref{prop.perm.extend.YtX.inj-sur} \textbf{(a)} (applied to $X\setminus
\left\{  x\right\}  $, $\gamma\circ\varepsilon$ and $\gamma^{-1}$ instead of
$Y$, $\alpha$ and $\beta$) yields%
\begin{align*}
\left(  \gamma\circ\varepsilon\circ\gamma^{-1}\right)  ^{\left(
X\setminus\left\{  x\right\}  \rightarrow X\right)  }  &  =\underbrace{\left(
\gamma\circ\varepsilon\right)  ^{\left(  X\setminus\left\{  x\right\}
\rightarrow X\right)  }}_{=\gamma^{\left(  X\setminus\left\{  x\right\}
\rightarrow X\right)  }\circ\varepsilon^{\left(  X\setminus\left\{  x\right\}
\rightarrow X\right)  }}\circ\underbrace{\left(  \gamma^{-1}\right)  ^{\left(
X\setminus\left\{  x\right\}  \rightarrow X\right)  }}_{=\left(
\gamma^{\left(  X\setminus\left\{  x\right\}  \rightarrow X\right)  }\right)
^{-1}}\\
&  =\gamma^{\left(  X\setminus\left\{  x\right\}  \rightarrow X\right)  }%
\circ\varepsilon^{\left(  X\setminus\left\{  x\right\}  \rightarrow X\right)
}\circ\left(  \gamma^{\left(  X\setminus\left\{  x\right\}  \rightarrow
X\right)  }\right)  ^{-1}.
\end{align*}
This proves Claim 4 \textbf{(c)}.]

We are now ready to prove Claim 2:

[\textit{Proof of Claim 2:} Since $x\in X$, we have $\left\vert X\setminus
\left\{  x\right\}  \right\vert =\left\vert X\right\vert -1<\left\vert
X\right\vert =N$.

Applying Proposition \ref{prop.perm.extend.YtX.inj-sur} \textbf{(d)} to
$Y=X\setminus\left\{  x\right\}  $, we obtain%
\begin{align}
&  \left\{  \delta^{\left(  X\setminus\left\{  x\right\}  \rightarrow
X\right)  }\ \mid\ \delta\in S_{X\setminus\left\{  x\right\}  }\right\}
\nonumber\\
&  =\left\{  \tau\in S_{X}\ \mid\ \tau\left(  z\right)  =z\text{ for every
}z\in X\setminus\left(  X\setminus\left\{  x\right\}  \right)  \right\}  .
\label{pf.prop.perm.extend.conj-inv.1.ext}%
\end{align}

We know that $\pi\in S_{X}$. In other words, $\pi$ is a permutation of $X$
(since $S_{X}$ is the set of all permutations of $X$). In other words, $\pi$
is a bijective map $X\rightarrow X$.

We are in one of the following two cases:

\textit{Case 1:} We have $\pi\left(  x\right)  =x$.

\textit{Case 2:} We don't have $\pi\left(  x\right)  =x$.

Let us first consider Case 1. In this case, we have $\pi\left(  x\right)  =x$.
Thus, $\pi\left(  z\right)  =z$ for every $z\in X\setminus\left(
X\setminus\left\{  x\right\}  \right)  $\ \ \ \ \footnote{\textit{Proof.} Let
$z\in X\setminus\left(  X\setminus\left\{  x\right\}  \right)  $. Then, $z\in
X\setminus\left(  X\setminus\left\{  x\right\}  \right)  \subseteq\left\{
x\right\}  $, so that $z=x$. Hence, $\pi\left(  z\right)  =\pi\left(
x\right)  =x=z$, qed.}. Hence, $\pi$ is a $\tau\in S_{X}$ satisfying
\newline$\left(  \tau\left(  z\right)  =z\text{ for every }z\in X\setminus
\left(  X\setminus\left\{  x\right\}  \right)  \right)  $. In other words,
\newline$\pi\in\left\{  \tau\in S_{X}\ \mid\ \tau\left(  z\right)  =z\text{
for every }z\in X\setminus\left(  X\setminus\left\{  x\right\}  \right)
\right\}  $. Hence,%
\begin{align*}
\pi &  \in\left\{  \tau\in S_{X}\ \mid\ \tau\left(  z\right)  =z\text{ for
every }z\in X\setminus\left(  X\setminus\left\{  x\right\}  \right)  \right\}
\\
&  =\left\{  \delta^{\left(  X\setminus\left\{  x\right\}  \rightarrow
X\right)  }\ \mid\ \delta\in S_{X\setminus\left\{  x\right\}  }\right\}
\ \ \ \ \ \ \ \ \ \ \left(  \text{by (\ref{pf.prop.perm.extend.conj-inv.1.ext}%
)}\right) \\
&  =\left\{  \left(  \pi^{\prime}\right)  ^{\left(  X\setminus\left\{
x\right\}  \rightarrow X\right)  }\ \mid\ \pi^{\prime}\in S_{X\setminus
\left\{  x\right\}  }\right\}
\end{align*}
(here, we have renamed the index $\delta$ as $\pi^{\prime}$). In other words,
there exists a $\pi^{\prime}\in S_{X\setminus\left\{  x\right\}  }$ such that
$\pi=\left(  \pi^{\prime}\right)  ^{\left(  X\setminus\left\{  x\right\}
\rightarrow X\right)  }$. Consider this $\pi^{\prime}$.

Proposition \ref{prop.perm.extend.YtX.inj-sur} \textbf{(c)} (applied to
$Y=X\setminus\left\{  x\right\}  $ and $\sigma=\pi^{\prime}$) shows that
$\pi^{\prime}$ satisfies $\left(  \pi^{\prime}\right)  ^{\left(
X\setminus\left\{  x\right\}  \rightarrow X\right)  }\in S_{X}$ and%
\[
\left(  \left(  \pi^{\prime}\right)  ^{-1}\right)  ^{\left(  X\setminus
\left\{  x\right\}  \rightarrow X\right)  }=\left(  \underbrace{\left(
\pi^{\prime}\right)  ^{\left(  X\setminus\left\{  x\right\}  \rightarrow
X\right)  }}_{=\pi}\right)  ^{-1}=\pi^{-1}.
\]

Now, Claim 3 (applied to $X\setminus\left\{  x\right\}  $ and $\pi^{\prime}$
instead of $X$ and $\pi$) yields that there exists a $\sigma\in S_{X\setminus
\left\{  x\right\}  }$ such that $\sigma\circ\pi^{\prime}\circ\sigma
^{-1}=\left(  \pi^{\prime}\right)  ^{-1}$. Denote this $\sigma$ by
$\sigma^{\prime}$. Thus, $\sigma^{\prime}\in S_{X\setminus\left\{  x\right\}
}$ and $\sigma^{\prime}\circ\pi^{\prime}\circ\left(  \sigma^{\prime}\right)
^{-1}=\left(  \pi^{\prime}\right)  ^{-1}$.

\begin{verlong}
We have $\sigma^{\prime}\in S_{X\setminus\left\{  x\right\}  }$. In other
words, $\sigma^{\prime}$ is a permutation of $X\setminus\left\{  x\right\}  $
(since $S_{X\setminus\left\{  x\right\}  }$ is the set of all permutations of
$X\setminus\left\{  x\right\}  $). In other words, $\sigma^{\prime}$ is a
bijective map $X\setminus\left\{  x\right\}  \rightarrow X\setminus\left\{
x\right\}  $.
\end{verlong}

Define the map $\sigma^{\prime\prime}:X\rightarrow X$ by
\[
\sigma^{\prime\prime}=\left(  \sigma^{\prime}\right)  ^{\left(  X\setminus
\left\{  x\right\}  \rightarrow X\right)  }.
\]
Thus, $\sigma^{\prime\prime}=\delta^{\left(  X\setminus\left\{  x\right\}
\rightarrow X\right)  }$ for some $\delta\in S_{X\setminus\left\{  x\right\}
}$ (namely, for $\delta=\sigma^{\prime}$). In other words,%
\begin{equation}
\sigma^{\prime\prime}\in\left\{  \delta^{\left(  X\setminus\left\{  x\right\}
\rightarrow X\right)  }\ \mid\ \delta\in S_{X\setminus\left\{  x\right\}
}\right\}  . \label{pf.prop.perm.extend.conj-inv.1.c1.4}%
\end{equation}

Now, Claim 4 \textbf{(c)} (applied to $\gamma=\sigma^{\prime}$ and
$\varepsilon=\pi^{\prime}$) yields%
\begin{align*}
&  \left(  \sigma^{\prime}\circ\pi^{\prime}\circ\left(  \sigma^{\prime
}\right)  ^{-1}\right)  ^{\left(  X\setminus\left\{  x\right\}  \rightarrow
X\right)  }\\
&  =\left(  \underbrace{\left(  \sigma^{\prime}\right)  ^{\left(
X\setminus\left\{  x\right\}  \rightarrow X\right)  }}_{=\sigma^{\prime\prime
}}\right)  \circ\left(  \underbrace{\left(  \pi^{\prime}\right)  ^{\left(
X\setminus\left\{  x\right\}  \rightarrow X\right)  }}_{=\pi}\right)
\circ\left(  \underbrace{\left(  \sigma^{\prime}\right)  ^{\left(
X\setminus\left\{  x\right\}  \rightarrow X\right)  }}_{=\sigma^{\prime\prime
}}\right)  ^{-1}\\
&  =\sigma^{\prime\prime}\circ\pi\circ\left(  \sigma^{\prime\prime}\right)
^{-1}.
\end{align*}
Compared with%
\[
\left(  \underbrace{\sigma^{\prime}\circ\pi^{\prime}\circ\left(
\sigma^{\prime}\right)  ^{-1}}_{=\left(  \pi^{\prime}\right)  ^{-1}}\right)
^{\left(  X\setminus\left\{  x\right\}  \rightarrow X\right)  }=\left(
\left(  \pi^{\prime}\right)  ^{-1}\right)  ^{\left(  X\setminus\left\{
x\right\}  \rightarrow X\right)  }=\pi^{-1},
\]
this yields $\sigma^{\prime\prime}\circ\pi\circ\left(  \sigma^{\prime\prime
}\right)  ^{-1}=\pi^{-1}$. So we have shown that $\sigma^{\prime\prime}$ is an
element of $\left\{  \delta^{\left(  X\setminus\left\{  x\right\}  \rightarrow
X\right)  }\ \mid\ \delta\in S_{X\setminus\left\{  x\right\}  }\right\}  $ (by
(\ref{pf.prop.perm.extend.conj-inv.1.c1.4})) and satisfies $\sigma
^{\prime\prime}\circ\pi\circ\left(  \sigma^{\prime\prime}\right)  ^{-1}%
=\pi^{-1}$. Hence, there exists a $\sigma\in\left\{  \delta^{\left(
X\setminus\left\{  x\right\}  \rightarrow X\right)  }\ \mid\ \delta\in
S_{X\setminus\left\{  x\right\}  }\right\}  $ such that $\sigma\circ\pi
\circ\sigma^{-1}=\pi^{-1}$ (namely, $\sigma=\sigma^{\prime\prime}$). Thus,
Claim 2 holds in Case 1.

Now, let us consider Case 2. In this case, we don't have $\pi\left(  x\right)
=x$. Hence, we have $\pi\left(  x\right)  \neq x$.

\begin{vershort}
Let $y=\pi\left(  x\right)  $. Combining $y=\pi\left(  x\right)  \in X$ with
$y=\pi\left(  x\right)  \neq x$, we obtain $y\in X\setminus\left\{  x\right\}
$. Hence, $x$ and $y$ are two distinct elements of $X$. Thus, the
transposition $t_{x,y}$ of $X$ is well-defined (according to Definition
\ref{def.transposX}). Consider this transposition $t_{x,y}$.
\end{vershort}

\begin{verlong}
Combining $\pi\left(  x\right)  \in X$ with $\pi\left(  x\right)  \neq x$, we
obtain $\pi\left(  x\right)  \in X\setminus\left\{  x\right\}  $. We can thus
define an element $y\in X\setminus\left\{  x\right\}  $ by $y=\pi\left(
x\right)  $. Consider this $y$. We have $y\in X\setminus\left\{  x\right\}
\subseteq X$; that is, $y$ is an element of $X$. Also, $y=\pi\left(  x\right)
\neq x$; thus, the elements $x$ and $y$ of $X$ are distinct. Hence, the
transposition $t_{x,y}$ of $X$ is well-defined (according to Definition
\ref{def.transposX}). Consider this transposition $t_{x,y}$. Note that
$\left\{  x,y\right\}  $ is a subset of $X$ (since $x$ and $y$ are elements of
$X$).
\end{verlong}

Lemma \ref{lem.sol.perm.transX.tij1} \textbf{(a)} (applied to $i=x$ and $j=y$)
yields $t_{x,y}\left(  x\right)  =y$. Lemma \ref{lem.sol.perm.transX.tij1}
\textbf{(b)} (applied to $i=x$ and $j=y$) yields $t_{x,y}\left(  y\right)
=x$. Lemma \ref{lem.sol.perm.transX.tij1} \textbf{(d)} (applied to $i=x$ and
$j=y$) yields $t_{x,y}\circ t_{x,y}=\operatorname*{id}$. Thus, $\left(
t_{x,y}\right)  ^{-1}=t_{x,y}$. Finally, Lemma \ref{lem.sol.perm.transX.tij1}
\textbf{(c)} (applied to $i=x$ and $j=y$) shows that we have%
\begin{equation}
t_{x,y}\left(  k\right)  =k\text{ for each }k\in X\setminus\left\{
x,y\right\}  . \label{pf.prop.perm.extend.conj-inv.1.txyk}%
\end{equation}
Thus,%
\begin{equation}
t_{x,y}\in\left\{  \alpha^{\left(  \left\{  x,y\right\}  \rightarrow X\right)
}\ \mid\ \alpha\in S_{\left\{  x,y\right\}  }\right\}
\label{pf.prop.perm.extend.conj-inv.1.txy}%
\end{equation}
\footnote{\textit{Proof of (\ref{pf.prop.perm.extend.conj-inv.1.txy}):} We
know that $t_{x,y}$ is a permutation of $X$. Hence, $t_{x,y}$ belongs to
$S_{X}$ (since $S_{X}$ is the set of all permutations of $X$) and satisfies
$t_{x,y}\left(  z\right)  =z$ for every $z\in X\setminus\left\{  x,y\right\}
$ (by (\ref{pf.prop.perm.extend.conj-inv.1.txyk}), applied to $k=z$). In other
words, $t_{x,y}$ is a $\tau\in S_{X}$ satisfying $\left(  \tau\left(
z\right)  =z\text{ for every }z\in X\setminus\left\{  x,y\right\}  \right)  $.
In other words, $t_{x,y}\in\left\{  \tau\in S_{X}\ \mid\ \tau\left(  z\right)
=z\text{ for every }z\in X\setminus\left\{  x,y\right\}  \right\}  $. But
Proposition \ref{prop.perm.extend.YtX.inj-sur} \textbf{(d)} (applied to
$\left\{  x,y\right\}  $ instead of $Y$) yields%
\[
\left\{  \delta^{\left(  \left\{  x,y\right\}  \rightarrow X\right)  }%
\ \mid\ \delta\in S_{\left\{  x,y\right\}  }\right\}  =\left\{  \tau\in
S_{X}\ \mid\ \tau\left(  z\right)  =z\text{ for every }z\in X\setminus\left\{
x,y\right\}  \right\}  .
\]
Now,%
\begin{align*}
t_{x,y}  &  \in\left\{  \tau\in S_{X}\ \mid\ \tau\left(  z\right)  =z\text{
for every }z\in X\setminus\left\{  x,y\right\}  \right\} \\
&  =\left\{  \delta^{\left(  \left\{  x,y\right\}  \rightarrow X\right)
}\ \mid\ \delta\in S_{\left\{  x,y\right\}  }\right\}  =\left\{
\alpha^{\left(  \left\{  x,y\right\}  \rightarrow X\right)  }\ \mid\ \alpha\in
S_{\left\{  x,y\right\}  }\right\}
\end{align*}
(here, we have renamed the index $\delta$ as $\alpha$). This proves
(\ref{pf.prop.perm.extend.conj-inv.1.txy}).}. In other words, there exists
some $\alpha\in S_{\left\{  x,y\right\}  }$ such that
\begin{equation}
t_{x,y}=\alpha^{\left(  \left\{  x,y\right\}  \rightarrow X\right)  }.
\label{pf.prop.perm.extend.conj-inv.1.txy=}%
\end{equation}
Consider this $\alpha$. (Of course, we can easily tell what this $\alpha$ is:
It is the permutation of $\left\{  x,y\right\}  $ that swaps $x$ with $y$. But
we don't need to know this.)

\begin{verlong}
We have $\alpha\in S_{\left\{  x,y\right\}  }$. In other words, $\alpha$ is a
permutation of the set $\left\{  x,y\right\}  $ (since $S_{\left\{
x,y\right\}  }$ is the set of all permutations of the set $\left\{
x,y\right\}  $). In other words, $\alpha$ is a bijective map from $\left\{
x,y\right\}  $ to $\left\{  x,y\right\}  $.
\end{verlong}

\begin{vershort}
Both $t_{x,y}$ and $\pi$ are permutations of $X$. Thus, their composition
$t_{x,y}\circ\pi$ is a permutation of $X$ as well. In other words,
$t_{x,y}\circ\pi\in S_{X}$.
\end{vershort}

\begin{verlong}
Both $t_{x,y}$ and $\pi$ are permutations of $X$. Thus, their composition
$t_{x,y}\circ\pi$ is a permutation of $X$ as well. In other words,
$t_{x,y}\circ\pi\in S_{X}$ (since $S_{X}$ is the set of all permutations of
$X$).
\end{verlong}

Set $\widetilde{\varepsilon}=t_{x,y}\circ\pi$. Then, $\widetilde{\varepsilon
}=t_{x,y}\circ\pi\in S_{X}$. Also,%
\[
\underbrace{\widetilde{\varepsilon}}_{=t_{x,y}\circ\pi}\left(  x\right)
=\left(  t_{x,y}\circ\pi\right)  \left(  x\right)  =t_{x,y}\left(
\underbrace{\pi\left(  x\right)  }_{=y}\right)  =t_{x,y}\left(  y\right)  =x.
\]
Thus, $\widetilde{\varepsilon}\left(  z\right)  =z$ for every $z\in
X\setminus\left(  X\setminus\left\{  x\right\}  \right)  $%
\ \ \ \ \footnote{\textit{Proof.} Let $z\in X\setminus\left(  X\setminus
\left\{  x\right\}  \right)  $. Then, $z\in X\setminus\left(  X\setminus
\left\{  x\right\}  \right)  \subseteq\left\{  x\right\}  $, so that $z=x$.
Hence, $\widetilde{\varepsilon}\left(  z\right)  =\widetilde{\varepsilon
}\left(  x\right)  =x=z$, qed.}. Hence, $\widetilde{\varepsilon}$ is a
$\tau\in S_{X}$ satisfying $\left(  \tau\left(  z\right)  =z\text{ for every
}z\in X\setminus\left(  X\setminus\left\{  x\right\}  \right)  \right)  $.
Thus,%
\begin{align*}
\widetilde{\varepsilon}  &  \in\left\{  \tau\in S_{X}\ \mid\ \tau\left(
z\right)  =z\text{ for every }z\in X\setminus\left(  X\setminus\left\{
x\right\}  \right)  \right\} \\
&  =\left\{  \delta^{\left(  X\setminus\left\{  x\right\}  \rightarrow
X\right)  }\ \mid\ \delta\in S_{X\setminus\left\{  x\right\}  }\right\}
\ \ \ \ \ \ \ \ \ \ \left(  \text{by (\ref{pf.prop.perm.extend.conj-inv.1.ext}%
)}\right) \\
&  =\left\{  \varepsilon^{\left(  X\setminus\left\{  x\right\}  \rightarrow
X\right)  }\ \mid\ \varepsilon\in S_{X\setminus\left\{  x\right\}  }\right\}
\end{align*}
(here, we have renamed the index $\delta$ as $\varepsilon$). In other words,
there exists some $\varepsilon\in S_{X\setminus\left\{  x\right\}  }$ such
that%
\begin{equation}
\widetilde{\varepsilon}=\varepsilon^{\left(  X\setminus\left\{  x\right\}
\rightarrow X\right)  }. \label{pf.prop.perm.extend.conj-inv.1.epsti}%
\end{equation}
Consider this $\varepsilon$.

Proposition \ref{prop.perm.extend.YtX.inj-sur} \textbf{(c)} (applied to
$Y=X\setminus\left\{  x\right\}  $ and $\sigma=\varepsilon$) shows that
$\varepsilon$ satisfies $\varepsilon^{\left(  X\setminus\left\{  x\right\}
\rightarrow X\right)  }\in S_{X}$ and
\begin{equation}
\left(  \varepsilon^{-1}\right)  ^{\left(  X\setminus\left\{  x\right\}
\rightarrow X\right)  }=\left(  \underbrace{\varepsilon^{\left(
X\setminus\left\{  x\right\}  \rightarrow X\right)  }}%
_{\substack{=\widetilde{\varepsilon}\\\text{(by
(\ref{pf.prop.perm.extend.conj-inv.1.epsti}))}}}\right)  ^{-1}%
=\widetilde{\varepsilon}^{-1}. \label{pf.prop.perm.extend.conj-inv.1.epsi-1}%
\end{equation}

Now, Claim 1 (applied to $X\setminus\left\{  x\right\}  $, $\varepsilon$ and
$y$ instead of $X$, $\pi$ and $x$) yields that there exists a $\sigma
\in\left\{  \delta^{\left(  \left(  X\setminus\left\{  x\right\}  \right)
\setminus\left\{  y\right\}  \rightarrow X\setminus\left\{  x\right\}
\right)  }\ \mid\ \delta\in S_{\left(  X\setminus\left\{  x\right\}  \right)
\setminus\left\{  y\right\}  }\right\}  $ such that $\sigma\circ
\varepsilon\circ\sigma^{-1}=\varepsilon^{-1}$. Denote this $\sigma$ by
$\gamma$. Then,%
\[
\gamma\in\left\{  \delta^{\left(  \left(  X\setminus\left\{  x\right\}
\right)  \setminus\left\{  y\right\}  \rightarrow X\setminus\left\{
x\right\}  \right)  }\ \mid\ \delta\in S_{\left(  X\setminus\left\{
x\right\}  \right)  \setminus\left\{  y\right\}  }\right\}  \text{ and }%
\gamma\circ\varepsilon\circ\gamma^{-1}=\varepsilon^{-1}.
\]

We have%
\begin{align*}
\gamma &  \in\left\{  \delta^{\left(  \left(  X\setminus\left\{  x\right\}
\right)  \setminus\left\{  y\right\}  \rightarrow X\setminus\left\{
x\right\}  \right)  }\ \mid\ \delta\in S_{\left(  X\setminus\left\{
x\right\}  \right)  \setminus\left\{  y\right\}  }\right\} \\
&  =\left\{  \delta^{\left(  X\setminus\left\{  x,y\right\}  \rightarrow
X\setminus\left\{  x\right\}  \right)  }\ \mid\ \delta\in S_{X\setminus
\left\{  x,y\right\}  }\right\}  \ \ \ \ \ \ \ \ \ \ \left(  \text{since
}\left(  X\setminus\left\{  x\right\}  \right)  \setminus\left\{  y\right\}
=X\setminus\left\{  x,y\right\}  \right) \\
&  =\left\{  \beta^{\left(  X\setminus\left\{  x,y\right\}  \rightarrow
X\setminus\left\{  x\right\}  \right)  }\ \mid\ \beta\in S_{X\setminus\left\{
x,y\right\}  }\right\}
\end{align*}
(here, we have renamed the index $\delta$ as $\beta$). In other words, there
exists some $\beta\in S_{X\setminus\left\{  x,y\right\}  }$ such that
$\gamma=\beta^{\left(  X\setminus\left\{  x,y\right\}  \rightarrow
X\setminus\left\{  x\right\}  \right)  }$. Consider this $\beta$.

Clearly, $X\setminus\left\{  x,y\right\}  $ is a subset of $X\setminus\left\{
x\right\}  $ (since $X\setminus\left\{  x,y\right\}  =\left(  X\setminus
\left\{  x\right\}  \right)  \setminus\left\{  y\right\}  $). Hence, from
Proposition \ref{prop.perm.extend.YtX.inj-sur} \textbf{(c)} (applied to
$X\setminus\left\{  x\right\}  $, $X\setminus\left\{  x,y\right\}  $ and
$\beta$ instead of $X$, $Y$ and $\sigma$), we obtain that $\beta$ satisfies
$\beta^{\left(  X\setminus\left\{  x,y\right\}  \rightarrow X\setminus\left\{
x\right\}  \right)  }\in S_{X\setminus\left\{  x\right\}  }$ and $\left(
\beta^{-1}\right)  ^{\left(  X\setminus\left\{  x,y\right\}  \rightarrow
X\setminus\left\{  x\right\}  \right)  }=\left(  \beta^{\left(  X\setminus
\left\{  x,y\right\}  \rightarrow X\setminus\left\{  x\right\}  \right)
}\right)  ^{-1}$. Thus,%
\[
\gamma=\beta^{\left(  X\setminus\left\{  x,y\right\}  \rightarrow
X\setminus\left\{  x\right\}  \right)  }\in S_{X\setminus\left\{  x\right\}
}.
\]

\begin{verlong}
Note that $\beta\in S_{X\setminus\left\{  x,y\right\}  }$. In other words,
$\beta$ is a permutation of $X\setminus\left\{  x,y\right\}  $ (since
$S_{X\setminus\left\{  x,y\right\}  }$ is the set of all permutations of
$X\setminus\left\{  x,y\right\}  $). In other words, $\beta$ is a bijective
map $X\setminus\left\{  x,y\right\}  \rightarrow X\setminus\left\{
x,y\right\}  $.
\end{verlong}

Recall that $\gamma\in S_{X\setminus\left\{  x\right\}  }$. Thus, Proposition
\ref{prop.perm.extend.YtX.inj-sur} \textbf{(c)} (applied to $Y=X\setminus
\left\{  x\right\}  $ and $\sigma=\gamma$) yields that $\gamma$ satisfies
$\gamma^{\left(  X\setminus\left\{  x\right\}  \rightarrow X\right)  }\in
S_{X}$ and $\left(  \gamma^{-1}\right)  ^{\left(  X\setminus\left\{
x\right\}  \rightarrow X\right)  }=\left(  \gamma^{\left(  X\setminus\left\{
x\right\}  \rightarrow X\right)  }\right)  ^{-1}$.

Now, set
\begin{equation}
\widetilde{\gamma}=\gamma^{\left(  X\setminus\left\{  x\right\}  \rightarrow
X\right)  }. \label{pf.prop.perm.extend.conj-inv.1.gamti=}%
\end{equation}
Thus,%
\begin{align}
\widetilde{\gamma}  &  =\gamma^{\left(  X\setminus\left\{  x\right\}
\rightarrow X\right)  }=\left(  \beta^{\left(  X\setminus\left\{  x,y\right\}
\rightarrow X\setminus\left\{  x\right\}  \right)  }\right)  ^{\left(
X\setminus\left\{  x\right\}  \rightarrow X\right)  }%
\ \ \ \ \ \ \ \ \ \ \left(  \text{since }\gamma=\beta^{\left(  X\setminus
\left\{  x,y\right\}  \rightarrow X\setminus\left\{  x\right\}  \right)
}\right) \nonumber\\
&  =\beta^{\left(  X\setminus\left\{  x,y\right\}  \rightarrow X\right)  }
\label{pf.prop.perm.extend.conj-inv.1.gamti}%
\end{align}
(by Proposition \ref{prop.perm.extend.YtX.ZtY}, applied to $Y=X\setminus
\left\{  x\right\}  $ and $Z=X\setminus\left\{  x,y\right\}  $ and
$\sigma=\beta$).

Now, Proposition \ref{prop.perm.extend.YtX.commute} (applied to $Y=\left\{
x,y\right\}  $) yields%
\[
\alpha^{\left(  \left\{  x,y\right\}  \rightarrow X\right)  }\circ
\beta^{\left(  X\setminus\left\{  x,y\right\}  \rightarrow X\right)  }%
=\beta^{\left(  X\setminus\left\{  x,y\right\}  \rightarrow X\right)  }%
\circ\alpha^{\left(  \left\{  x,y\right\}  \rightarrow X\right)  }.
\]
In view of (\ref{pf.prop.perm.extend.conj-inv.1.txy=}) and
(\ref{pf.prop.perm.extend.conj-inv.1.gamti}), this rewrites as
\begin{equation}
t_{x,y}\circ\widetilde{\gamma}=\widetilde{\gamma}\circ t_{x,y}.
\label{pf.prop.perm.extend.conj-inv.1.commutat}%
\end{equation}

\begin{vershort}
Both $\gamma$ and $\varepsilon$ are elements of $S_{X\setminus\left\{
x\right\}  }$. In other words, both $\gamma$ and $\varepsilon$ are
permutations of $X\setminus\left\{  x\right\}  $. Hence, $\gamma
\circ\varepsilon$ is a permutation of $X\setminus\left\{  x\right\}  $ as
well. In other words, $\gamma\circ\varepsilon\in S_{X\setminus\left\{
x\right\}  }$.
\end{vershort}

\begin{verlong}
Recall that $\gamma\in S_{X\setminus\left\{  x\right\}  }$. Thus, both
$\gamma$ and $\varepsilon$ are elements of $S_{X\setminus\left\{  x\right\}
}$. In other words, both $\gamma$ and $\varepsilon$ are permutations of
$X\setminus\left\{  x\right\}  $ (since $S_{X\setminus\left\{  x\right\}  }$
is the set of all permutations of $X\setminus\left\{  x\right\}  $). In other
words, both $\gamma$ and $\varepsilon$ are bijective maps $X\setminus\left\{
x\right\}  \rightarrow X\setminus\left\{  x\right\}  $. Hence, $\gamma
\circ\varepsilon$ is a bijective map $X\setminus\left\{  x\right\}
\rightarrow X\setminus\left\{  x\right\}  $ as well, i.e., is a permutation of
$X\setminus\left\{  x\right\}  $ as well, i.e., is an element of
$S_{X\setminus\left\{  x\right\}  }$. In other words, $\gamma\circ
\varepsilon\in S_{X\setminus\left\{  x\right\}  }$.
\end{verlong}

Moreover, Claim 4 \textbf{(b)} yields%
\[
\left(  \gamma\circ\varepsilon\right)  ^{\left(  X\setminus\left\{  x\right\}
\rightarrow X\right)  }=\underbrace{\gamma^{\left(  X\setminus\left\{
x\right\}  \rightarrow X\right)  }}_{\substack{=\widetilde{\gamma}\\\text{(by
(\ref{pf.prop.perm.extend.conj-inv.1.gamti=}))}}}\circ\underbrace{\varepsilon
^{\left(  X\setminus\left\{  x\right\}  \rightarrow X\right)  }}%
_{\substack{=\widetilde{\varepsilon}\\\text{(by
(\ref{pf.prop.perm.extend.conj-inv.1.epsti}))}}}=\widetilde{\gamma}%
\circ\widetilde{\varepsilon},
\]
so that $\widetilde{\gamma}\circ\widetilde{\varepsilon}=\left(  \gamma
\circ\varepsilon\right)  ^{\left(  X\setminus\left\{  x\right\}  \rightarrow
X\right)  }$. Hence, $\widetilde{\gamma}\circ\widetilde{\varepsilon}%
=\delta^{\left(  X\setminus\left\{  x\right\}  \rightarrow X\right)  }$ for
some $\delta\in S_{X\setminus\left\{  x\right\}  }$ (namely, $\delta
=\gamma\circ\varepsilon$), because $\gamma\circ\varepsilon\in S_{X\setminus
\left\{  x\right\}  }$. In other words,%
\begin{equation}
\widetilde{\gamma}\circ\widetilde{\varepsilon}\in\left\{  \delta^{\left(
X\setminus\left\{  x\right\}  \rightarrow X\right)  }\ \mid\ \delta\in
S_{X\setminus\left\{  x\right\}  }\right\}  .
\label{pf.prop.perm.extend.conj-inv.1.right-type}%
\end{equation}

On the other hand, recall that $\gamma\circ\varepsilon\circ\gamma
^{-1}=\varepsilon^{-1}$. Thus,
\[
\left(  \gamma\circ\varepsilon\circ\gamma^{-1}\right)  ^{\left(
X\setminus\left\{  x\right\}  \rightarrow X\right)  }=\left(  \varepsilon
^{-1}\right)  ^{\left(  X\setminus\left\{  x\right\}  \rightarrow X\right)
}=\widetilde{\varepsilon}^{-1}\ \ \ \ \ \ \ \ \ \ \left(  \text{by
(\ref{pf.prop.perm.extend.conj-inv.1.epsi-1})}\right)  .
\]
Compared with%
\begin{align*}
&  \left(  \gamma\circ\varepsilon\circ\gamma^{-1}\right)  ^{\left(
X\setminus\left\{  x\right\}  \rightarrow X\right)  }\\
&  =\left(  \underbrace{\gamma^{\left(  X\setminus\left\{  x\right\}
\rightarrow X\right)  }}_{\substack{=\widetilde{\gamma}\\\text{(by
(\ref{pf.prop.perm.extend.conj-inv.1.gamti=}))}}}\right)  \circ\left(
\underbrace{\varepsilon^{\left(  X\setminus\left\{  x\right\}  \rightarrow
X\right)  }}_{\substack{=\widetilde{\varepsilon}\\\text{(by
(\ref{pf.prop.perm.extend.conj-inv.1.epsti}))}}}\right)  \circ\left(
\underbrace{\gamma^{\left(  X\setminus\left\{  x\right\}  \rightarrow
X\right)  }}_{\substack{=\widetilde{\gamma}\\\text{(by
(\ref{pf.prop.perm.extend.conj-inv.1.gamti=}))}}}\right)  ^{-1}%
\ \ \ \ \ \ \ \ \ \ \left(  \text{by Claim 4 \textbf{(c)}}\right) \\
&  =\widetilde{\gamma}\circ\widetilde{\varepsilon}\circ\widetilde{\gamma}%
^{-1},
\end{align*}
this yields%
\begin{equation}
\widetilde{\varepsilon}^{-1}=\widetilde{\gamma}\circ\widetilde{\varepsilon
}\circ\widetilde{\gamma}^{-1}.
\label{pf.prop.perm.extend.conj-inv.1.conj-for-epsi}%
\end{equation}

But
\begin{align*}
&  \left(  \underbrace{t_{x,y}\circ\widetilde{\gamma}}%
_{\substack{=\widetilde{\gamma}\circ t_{x,y}\\\text{(by
(\ref{pf.prop.perm.extend.conj-inv.1.commutat}))}}}\right)  ^{-1}\circ\left(
\widetilde{\gamma}\circ\underbrace{\widetilde{\varepsilon}}_{=t_{x,y}\circ\pi
}\right) \\
&  =\left(  \widetilde{\gamma}\circ t_{x,y}\right)  ^{-1}\circ\left(
\widetilde{\gamma}\circ t_{x,y}\circ\pi\right)  =\underbrace{\left(
\widetilde{\gamma}\circ t_{x,y}\right)  ^{-1}\circ\left(  \widetilde{\gamma
}\circ t_{x,y}\right)  }_{=\operatorname*{id}}\circ\pi=\pi.
\end{align*}
Thus, $\pi=\left(  t_{x,y}\circ\widetilde{\gamma}\right)  ^{-1}\circ\left(
\widetilde{\gamma}\circ\widetilde{\varepsilon}\right)  $, so that%
\begin{align}
&  \left(  \widetilde{\gamma}\circ\widetilde{\varepsilon}\right)
\circ\underbrace{\pi}_{=\left(  t_{x,y}\circ\widetilde{\gamma}\right)
^{-1}\circ\left(  \widetilde{\gamma}\circ\widetilde{\varepsilon}\right)
}\circ\left(  \widetilde{\gamma}\circ\widetilde{\varepsilon}\right)
^{-1}\nonumber\\
&  =\left(  \widetilde{\gamma}\circ\widetilde{\varepsilon}\right)
\circ\left(  t_{x,y}\circ\widetilde{\gamma}\right)  ^{-1}\circ
\underbrace{\left(  \widetilde{\gamma}\circ\widetilde{\varepsilon}\right)
\circ\left(  \widetilde{\gamma}\circ\widetilde{\varepsilon}\right)  ^{-1}%
}_{=\operatorname*{id}}=\left(  \widetilde{\gamma}\circ\widetilde{\varepsilon
}\right)  \circ\underbrace{\left(  t_{x,y}\circ\widetilde{\gamma}\right)
^{-1}}_{=\widetilde{\gamma}^{-1}\circ\left(  t_{x,y}\right)  ^{-1}}\nonumber\\
&  =\underbrace{\left(  \widetilde{\gamma}\circ\widetilde{\varepsilon}\right)
\circ\widetilde{\gamma}^{-1}}_{\substack{=\widetilde{\gamma}\circ
\widetilde{\varepsilon}\circ\widetilde{\gamma}^{-1}=\widetilde{\varepsilon
}^{-1}\\\text{(by (\ref{pf.prop.perm.extend.conj-inv.1.conj-for-epsi}))}%
}}\circ\underbrace{\left(  t_{x,y}\right)  ^{-1}}_{=t_{x,y}}%
=\widetilde{\varepsilon}^{-1}\circ t_{x,y}.
\label{pf.prop.perm.extend.conj-inv.1.almost-there}%
\end{align}
Since $\widetilde{\varepsilon}=t_{x,y}\circ\pi$, we have
$\widetilde{\varepsilon}^{-1}=\left(  t_{x,y}\circ\pi\right)  ^{-1}=\pi
^{-1}\circ\left(  t_{x,y}\right)  ^{-1}$. Hence,
(\ref{pf.prop.perm.extend.conj-inv.1.almost-there}) becomes%
\[
\left(  \widetilde{\gamma}\circ\widetilde{\varepsilon}\right)  \circ\pi
\circ\left(  \widetilde{\gamma}\circ\widetilde{\varepsilon}\right)
^{-1}=\underbrace{\widetilde{\varepsilon}^{-1}}_{=\pi^{-1}\circ\left(
t_{x,y}\right)  ^{-1}}\circ t_{x,y}=\pi^{-1}\circ\underbrace{\left(
t_{x,y}\right)  ^{-1}\circ t_{x,y}}_{=\operatorname*{id}}=\pi^{-1}.
\]

Thus, $\widetilde{\gamma}\circ\widetilde{\varepsilon}$ is an element of
$\left\{  \delta^{\left(  X\setminus\left\{  x\right\}  \rightarrow X\right)
}\ \mid\ \delta\in S_{X\setminus\left\{  x\right\}  }\right\}  $ (by
(\ref{pf.prop.perm.extend.conj-inv.1.right-type})) and satisfies $\left(
\widetilde{\gamma}\circ\widetilde{\varepsilon}\right)  \circ\pi\circ\left(
\widetilde{\gamma}\circ\widetilde{\varepsilon}\right)  ^{-1}=\pi^{-1}$. Hence,
there exists a $\sigma\in\left\{  \delta^{\left(  X\setminus\left\{
x\right\}  \rightarrow X\right)  }\ \mid\ \delta\in S_{X\setminus\left\{
x\right\}  }\right\}  $ such that $\sigma\circ\pi\circ\sigma^{-1}=\pi^{-1}$
(namely, $\sigma=\widetilde{\gamma}\circ\widetilde{\varepsilon}$). Hence,
Claim 2 is proven in Case 2.

We have now proven that Claim 2 holds in each of the two Cases 1 and 2. Since
these two Cases cover all possibilities, we thus conclude that Claim 2 always holds.]

Thus, we have proven Claim 2. In other words, we have proven that Proposition
\ref{prop.perm.extend.conj-inv.1} holds whenever $\left\vert X\right\vert =N$.
This completes the induction step. The induction proof of Proposition
\ref{prop.perm.extend.conj-inv.1} is thus complete.
\end{proof}

\begin{proof}
[Proof of Theorem \ref{thm.perm.extend.conj-inv.2}.]If $X=\varnothing$, then
Theorem \ref{thm.perm.extend.conj-inv.2} holds for simple
reasons\footnote{\textit{Proof.} Assume that $X=\varnothing$. We must prove
that Theorem \ref{thm.perm.extend.conj-inv.2} holds.
\par
We have $X=\varnothing$. Thus, no $x\in X$ exists. Hence, $\left(  \pi\circ
\pi\circ\pi^{-1}\right)  \left(  x\right)  =\pi^{-1}\left(  x\right)  $ holds
for all $x\in X$ (indeed, this is vacuously true, since no $x\in X$ exists).
In other words, we have $\pi\circ\pi\circ\pi^{-1}=\pi^{-1}$. Thus, there
exists a $\sigma\in S_{X}$ such that $\sigma\circ\pi\circ\sigma^{-1}=\pi^{-1}$
(namely, $\sigma=\pi$). Hence, Theorem \ref{thm.perm.extend.conj-inv.2} holds.
\par
We thus have proven that if $X=\varnothing$, then Theorem
\ref{thm.perm.extend.conj-inv.2} holds. Qed.}. Hence, for the rest of this
proof of Theorem \ref{thm.perm.extend.conj-inv.2}, we can WLOG assume that we
don't have $X=\varnothing$. Assume this.

\begin{vershort}
There exists some $x\in X$ (since we don't have $X=\varnothing$). Consider
this $x$. Hence, Proposition \ref{prop.perm.extend.conj-inv.1} shows that
there exists a $\sigma\in\left\{  \delta^{\left(  X\setminus\left\{
x\right\}  \rightarrow X\right)  }\ \mid\ \delta\in S_{X\setminus\left\{
x\right\}  }\right\}  $ such that $\sigma\circ\pi\circ\sigma^{-1}=\pi^{-1}$.
Consider this $\sigma$, and denote it by $\alpha$. Thus, $\alpha$ is an
element of $\left\{  \delta^{\left(  X\setminus\left\{  x\right\}  \rightarrow
X\right)  }\ \mid\ \delta\in S_{X\setminus\left\{  x\right\}  }\right\}  $ and
satisfies $\alpha\circ\pi\circ\alpha^{-1}=\pi^{-1}$.
\end{vershort}

\begin{verlong}
We don't have $X=\varnothing$. Thus, the set $X$ is nonempty. Hence, there
exists some $x\in X$. Consider this $x$. Hence, Proposition
\ref{prop.perm.extend.conj-inv.1} shows that there exists a $\sigma\in\left\{
\delta^{\left(  X\setminus\left\{  x\right\}  \rightarrow X\right)  }%
\ \mid\ \delta\in S_{X\setminus\left\{  x\right\}  }\right\}  $ such that
$\sigma\circ\pi\circ\sigma^{-1}=\pi^{-1}$. Consider this $\sigma$, and denote
it by $\alpha$. Thus, $\alpha$ is an element of $\left\{  \delta^{\left(
X\setminus\left\{  x\right\}  \rightarrow X\right)  }\ \mid\ \delta\in
S_{X\setminus\left\{  x\right\}  }\right\}  $ and satisfies $\alpha\circ
\pi\circ\alpha^{-1}=\pi^{-1}$.
\end{verlong}

We know that $\alpha$ is an element of $\left\{  \delta^{\left(
X\setminus\left\{  x\right\}  \rightarrow X\right)  }\ \mid\ \delta\in
S_{X\setminus\left\{  x\right\}  }\right\}  $. In other words, there exists
some $\delta\in S_{X\setminus\left\{  x\right\}  }$ such that $\alpha
=\delta^{\left(  X\setminus\left\{  x\right\}  \rightarrow X\right)  }$.
Consider this $\delta$.

Now, Proposition \ref{prop.perm.extend.YtX.inj-sur} \textbf{(c)} (applied to
$Y=X\setminus\left\{  x\right\}  $ and $\sigma=\delta$) yields that $\delta$
satisfies $\delta^{\left(  X\setminus\left\{  x\right\}  \rightarrow X\right)
}\in S_{X}$ and $\left(  \delta^{-1}\right)  ^{\left(  X\setminus\left\{
x\right\}  \rightarrow X\right)  }=\left(  \delta^{\left(  X\setminus\left\{
x\right\}  \rightarrow X\right)  }\right)  ^{-1}$. Now, $\alpha=\delta
^{\left(  X\setminus\left\{  x\right\}  \rightarrow X\right)  }\in S_{X}$ and
$\alpha\circ\pi\circ\alpha^{-1}=\pi^{-1}$. Hence, there exists a $\sigma\in
S_{X}$ such that $\sigma\circ\pi\circ\sigma^{-1}=\pi^{-1}$ (namely,
$\sigma=\alpha$). This proves Theorem \ref{thm.perm.extend.conj-inv.2}.
\end{proof}

\begin{proof}
[Solution to Exercise \ref{exe.perm.extend.proofs}.]We have now proven
Proposition \ref{prop.perm.extend.YtX.wd}, Proposition
\ref{prop.perm.extend.YtX.inj-sur}, Proposition \ref{prop.perm.extend.YtX.ZtY}%
, Proposition \ref{prop.perm.extend.YtX.commute}, Proposition
\ref{prop.perm.extend.conj-inv.1} and Theorem \ref{thm.perm.extend.conj-inv.2}%
. Thus, Exercise \ref{exe.perm.extend.proofs} is solved.
\end{proof}

\subsection{Solution to Exercise \ref{exe.perm.aj-ai-sum}}

We shall use the Iverson bracket notation introduced in Definition
\ref{def.iverson}.

\begin{lemma}
\label{lem.sol.perm.aj-ai-sum.1}Let $n\in\mathbb{N}$. Let $\left[  n\right]
=\left\{  1,2,\ldots,n\right\}  $. Let $j\in\left[  n\right]  $.

\textbf{(a)} Then, $\sum_{i\in\left[  n\right]  }\left[  i\leq j\right]  =j$.

\textbf{(b)} Let $\sigma\in S_{n}$. Then, $\sum_{i\in\left[  n\right]
}\left[  \sigma\left(  i\right)  \leq\sigma\left(  j\right)  \right]
=\sigma\left(  j\right)  $.
\end{lemma}

\begin{vershort}
\begin{proof}
[Proof of Lemma \ref{lem.sol.perm.aj-ai-sum.1}.]\textbf{(a)} From $\left[
n\right]  =\left\{  1,2,\ldots,n\right\}  $, we obtain
\begin{align*}
\sum_{i\in\left[  n\right]  }\left[  i\leq j\right]   &  =\underbrace{\sum
_{i\in\left\{  1,2,\ldots,n\right\}  }}_{=\sum_{i=1}^{n}}\left[  i\leq
j\right]  =\sum_{i=1}^{n}\left[  i\leq j\right]  =\sum_{i=1}^{j}%
\underbrace{\left[  i\leq j\right]  }_{\substack{=1\\\text{(since we have
}i\leq j\text{)}}}+\sum_{i=j+1}^{n}\underbrace{\left[  i\leq j\right]
}_{\substack{=0\\\text{(since we don't have }i\leq j\\\text{(because }i\geq
j+1>j\text{))}}}\\
&  =\sum_{i=1}^{j}1+\underbrace{\sum_{i=j+1}^{n}0}_{=0}=\sum_{i=1}^{j}%
1=j\cdot1=j.
\end{align*}
This proves Lemma \ref{lem.sol.perm.aj-ai-sum.1} \textbf{(a)}.

\textbf{(b)} Recall that $S_{n}$ is the set of all permutations of the set
$\left\{  1,2,\ldots,n\right\}  $ (by the definition of $S_{n}$). In other
words, $S_{n}$ is the set of all permutations of the set $\left[  n\right]  $
(since $\left[  n\right]  =\left\{  1,2,\ldots,n\right\}  $). Thus, $\sigma$
is a permutation of the set $\left[  n\right]  $ (since $\sigma\in S_{n}$). In
other words, $\sigma$ is a bijection $\left[  n\right]  \rightarrow\left[
n\right]  $. Hence, we can substitute $i$ for $\sigma\left(  i\right)  $ in
the sum $\sum_{i\in\left[  n\right]  }\left[  \sigma\left(  i\right)
\leq\sigma\left(  j\right)  \right]  $. We thus obtain%
\[
\sum_{i\in\left[  n\right]  }\left[  \sigma\left(  i\right)  \leq\sigma\left(
j\right)  \right]  =\sum_{i\in\left[  n\right]  }\left[  i\leq\sigma\left(
j\right)  \right]  =\sigma\left(  j\right)
\]
(by Lemma \ref{lem.sol.perm.aj-ai-sum.1} \textbf{(a)}, applied to
$\sigma\left(  j\right)  $ instead of $j$). This proves Lemma
\ref{lem.sol.perm.aj-ai-sum.1} \textbf{(b)}.
\end{proof}
\end{vershort}

\begin{verlong}
\begin{proof}
[Proof of Lemma \ref{lem.sol.perm.aj-ai-sum.1}.]\textbf{(a)} We have
$j\in\left[  n\right]  =\left\{  1,2,\ldots,n\right\}  $, so that $1\leq j\leq
n$.

But $\left[  n\right]  =\left\{  1,2,\ldots,n\right\}  $. Thus,%
\begin{align*}
\sum_{i\in\left[  n\right]  }\left[  i\leq j\right]   &  =\underbrace{\sum
_{i\in\left\{  1,2,\ldots,n\right\}  }}_{=\sum_{i=1}^{n}}\left[  i\leq
j\right]  =\sum_{i=1}^{n}\left[  i\leq j\right]  =\sum_{i=1}^{j}%
\underbrace{\left[  i\leq j\right]  }_{\substack{=1\\\text{(since we have
}i\leq j\text{)}}}+\sum_{i=j+1}^{n}\underbrace{\left[  i\leq j\right]
}_{\substack{=0\\\text{(since we don't have }i\leq j\\\text{(because }i\geq
j+1>j\text{))}}}\\
&  \ \ \ \ \ \ \ \ \ \ \left(
\begin{array}
[c]{c}%
\text{here, we have split the summation at }i=j\text{,}\\
\text{since }1\leq j\leq n
\end{array}
\right) \\
&  =\sum_{i=1}^{j}1+\underbrace{\sum_{i=j+1}^{n}0}_{=0}=\sum_{i=1}^{j}%
1=j\cdot1=j.
\end{align*}
This proves Lemma \ref{lem.sol.perm.aj-ai-sum.1} \textbf{(a)}.

\textbf{(b)} Recall that $S_{n}$ is the set of all permutations of the set
$\left\{  1,2,\ldots,n\right\}  $ (by the definition of $S_{n}$). In other
words, $S_{n}$ is the set of all permutations of the set $\left[  n\right]  $
(since $\left[  n\right]  =\left\{  1,2,\ldots,n\right\}  $).

But $\sigma\in S_{n}$. In other words, $\sigma$ is a permutation of the set
$\left[  n\right]  $ (since $S_{n}$ is the set of all permutations of the set
$\left[  n\right]  $). In other words, $\sigma$ is a bijection $\left[
n\right]  \rightarrow\left[  n\right]  $. Hence, we can substitute $i$ for
$\sigma\left(  i\right)  $ in the sum $\sum_{i\in\left[  n\right]  }\left[
\sigma\left(  i\right)  \leq\sigma\left(  j\right)  \right]  $. We thus obtain%
\begin{equation}
\sum_{i\in\left[  n\right]  }\left[  \sigma\left(  i\right)  \leq\sigma\left(
j\right)  \right]  =\sum_{i\in\left[  n\right]  }\left[  i\leq\sigma\left(
j\right)  \right]  . \label{pf.lem.sol.perm.aj-ai-sum.1.b.1}%
\end{equation}

But $\sigma$ is a bijection $\left[  n\right]  \rightarrow\left[  n\right]  $.
Hence, $\sigma\left(  j\right)  \in\left[  n\right]  $ (since $j\in\left[
n\right]  $). Hence, Lemma \ref{lem.sol.perm.aj-ai-sum.1} \textbf{(a)}
(applied to $\sigma\left(  j\right)  $ instead of $j$) yields $\sum
_{i\in\left[  n\right]  }\left[  i\leq\sigma\left(  j\right)  \right]
=\sigma\left(  j\right)  $. Thus, (\ref{pf.lem.sol.perm.aj-ai-sum.1.b.1})
becomes
\[
\sum_{i\in\left[  n\right]  }\left[  \sigma\left(  i\right)  \leq\sigma\left(
j\right)  \right]  =\sum_{i\in\left[  n\right]  }\left[  i\leq\sigma\left(
j\right)  \right]  =\sigma\left(  j\right)  .
\]
This proves Lemma \ref{lem.sol.perm.aj-ai-sum.1} \textbf{(b)}.
\end{proof}
\end{verlong}

\begin{lemma}
\label{lem.sol.perm.aj-ai-sum.2}Let $n\in\mathbb{N}$. Let $\left[  n\right]
=\left\{  1,2,\ldots,n\right\}  $. Let $\sigma\in S_{n}$. Let $a_{1}%
,a_{2},\ldots,a_{n}$ be any $n$ numbers. (Here, \textquotedblleft
number\textquotedblright\ means \textquotedblleft real
number\textquotedblright\ or \textquotedblleft complex
number\textquotedblright\ or \textquotedblleft rational
number\textquotedblright, as you prefer; this makes no difference.) Then,%
\[
\sum_{\substack{1\leq i<j\leq n;\\\sigma\left(  i\right)  >\sigma\left(
j\right)  }}\left(  a_{j}-a_{i}\right)  =\sum_{i\in\left[  n\right]  }%
\sum_{j\in\left[  n\right]  }\left[  i<j\right]  \left[  \sigma\left(
i\right)  >\sigma\left(  j\right)  \right]  \left(  a_{j}-a_{i}\right)  .
\]
[Here, the summation sign \textquotedblleft$\sum_{\substack{1\leq i<j\leq
n;\\\sigma\left(  i\right)  >\sigma\left(  j\right)  }}$\textquotedblright%
\ means \textquotedblleft$\sum_{\substack{\left(  i,j\right)  \in\left\{
1,2,\ldots,n\right\}  ^{2};\\i<j\text{ and }\sigma\left(  i\right)
>\sigma\left(  j\right)  }}$\textquotedblright; this is a sum over all
inversions of $\sigma$.]
\end{lemma}

\begin{proof}
[Proof of Lemma \ref{lem.sol.perm.aj-ai-sum.2}.]If $\mathcal{A}$ and
$\mathcal{B}$ are two logical statements, then $\left[  \mathcal{A}%
\wedge\mathcal{B}\right]  =\left[  \mathcal{A}\right]  \left[  \mathcal{B}%
\right]  $ (by Exercise \ref{exe.iverson-prop} \textbf{(b)}). Hence, if
$\mathcal{A}$ and $\mathcal{B}$ are two logical statements, then%
\begin{equation}
\left[  \mathcal{A}\right]  \left[  \mathcal{B}\right]  =\left[
\mathcal{A}\wedge\mathcal{B}\right]  =\left[  \mathcal{A}\text{ and
}\mathcal{B}\right]  \label{pf.lem.sol.perm.aj-ai-sum.2.AandB}%
\end{equation}
(since \textquotedblleft$\mathcal{A}\wedge\mathcal{B}$\textquotedblright%
\ means \textquotedblleft$\mathcal{A}$ and $\mathcal{B}$\textquotedblright).

Recall that the summation sign \textquotedblleft$\sum_{\substack{1\leq i<j\leq
n;\\\sigma\left(  i\right)  >\sigma\left(  j\right)  }}$\textquotedblright%
\ means \textquotedblleft$\sum_{\substack{\left(  i,j\right)  \in\left\{
1,2,\ldots,n\right\}  ^{2};\\i<j\text{ and }\sigma\left(  i\right)
>\sigma\left(  j\right)  }}$\textquotedblright. Hence,%
\[
\sum_{\substack{1\leq i<j\leq n;\\\sigma\left(  i\right)  >\sigma\left(
j\right)  }}\left(  a_{j}-a_{i}\right)  =\sum_{\substack{\left(  i,j\right)
\in\left\{  1,2,\ldots,n\right\}  ^{2};\\i<j\text{ and }\sigma\left(
i\right)  >\sigma\left(  j\right)  }}\left(  a_{j}-a_{i}\right)
=\sum_{\substack{\left(  i,j\right)  \in\left[  n\right]  ^{2};\\i<j\text{ and
}\sigma\left(  i\right)  >\sigma\left(  j\right)  }}\left(  a_{j}%
-a_{i}\right)
\]
(since $\left\{  1,2,\ldots,n\right\}  =\left[  n\right]  $). Comparing this
with%
\begin{align*}
&  \underbrace{\sum_{i\in\left[  n\right]  }\sum_{j\in\left[  n\right]  }%
}_{=\sum_{\left(  i,j\right)  \in\left[  n\right]  ^{2}}}\underbrace{\left[
i<j\right]  \left[  \sigma\left(  i\right)  >\sigma\left(  j\right)  \right]
}_{\substack{=\left[  i<j\text{ and }\sigma\left(  i\right)  >\sigma\left(
j\right)  \right]  \\\text{(by (\ref{pf.lem.sol.perm.aj-ai-sum.2.AandB})
(applied}\\\text{to }\mathcal{A}=\left(  i<j\right)  \text{ and }%
\mathcal{B}=\left(  \sigma\left(  i\right)  >\sigma\left(  j\right)  \right)
\text{))}}}\left(  a_{j}-a_{i}\right) \\
&  =\sum_{\left(  i,j\right)  \in\left[  n\right]  ^{2}}\left[  i<j\text{ and
}\sigma\left(  i\right)  >\sigma\left(  j\right)  \right]  \left(  a_{j}%
-a_{i}\right) \\
&  =\sum_{\substack{\left(  i,j\right)  \in\left[  n\right]  ^{2};\\i<j\text{
and }\sigma\left(  i\right)  >\sigma\left(  j\right)  }}\underbrace{\left[
i<j\text{ and }\sigma\left(  i\right)  >\sigma\left(  j\right)  \right]
}_{\substack{=1\\\text{(since }i<j\text{ and }\sigma\left(  i\right)
>\sigma\left(  j\right)  \text{)}}}\left(  a_{j}-a_{i}\right) \\
&  \ \ \ \ \ \ \ \ \ \ +\sum_{\substack{\left(  i,j\right)  \in\left[
n\right]  ^{2};\\\text{not }\left(  i<j\text{ and }\sigma\left(  i\right)
>\sigma\left(  j\right)  \right)  }}\underbrace{\left[  i<j\text{ and }%
\sigma\left(  i\right)  >\sigma\left(  j\right)  \right]  }%
_{\substack{=0\\\text{(since we don't have }\left(  i<j\text{ and }%
\sigma\left(  i\right)  >\sigma\left(  j\right)  \right)  \text{)}}}\left(
a_{j}-a_{i}\right) \\
&  =\sum_{\substack{\left(  i,j\right)  \in\left[  n\right]  ^{2};\\i<j\text{
and }\sigma\left(  i\right)  >\sigma\left(  j\right)  }}\left(  a_{j}%
-a_{i}\right)  +\underbrace{\sum_{\substack{\left(  i,j\right)  \in\left[
n\right]  ^{2};\\\text{not }\left(  i<j\text{ and }\sigma\left(  i\right)
>\sigma\left(  j\right)  \right)  }}0\left(  a_{j}-a_{i}\right)  }_{=0}%
=\sum_{\substack{\left(  i,j\right)  \in\left[  n\right]  ^{2};\\i<j\text{ and
}\sigma\left(  i\right)  >\sigma\left(  j\right)  }}\left(  a_{j}%
-a_{i}\right)  ,
\end{align*}
we obtain
\[
\sum_{\substack{1\leq i<j\leq n;\\\sigma\left(  i\right)  >\sigma\left(
j\right)  }}\left(  a_{j}-a_{i}\right)  =\sum_{i\in\left[  n\right]  }%
\sum_{j\in\left[  n\right]  }\left[  i<j\right]  \left[  \sigma\left(
i\right)  >\sigma\left(  j\right)  \right]  \left(  a_{j}-a_{i}\right)  .
\]
This proves Lemma \ref{lem.sol.perm.aj-ai-sum.2}.
\end{proof}

\begin{lemma}
\label{lem.sol.perm.aj-ai-sum.3}Let $n\in\mathbb{N}$. Let $\left[  n\right]
=\left\{  1,2,\ldots,n\right\}  $. Let $\sigma\in S_{n}$. Let $i\in\left[
n\right]  $ and $j\in\left[  n\right]  $. Then,%
\[
\left[  i<j\right]  \left[  \sigma\left(  i\right)  >\sigma\left(  j\right)
\right]  -\left[  j<i\right]  \left[  \sigma\left(  j\right)  >\sigma\left(
i\right)  \right]  =\left[  i\leq j\right]  -\left[  \sigma\left(  i\right)
\leq\sigma\left(  j\right)  \right]  .
\]

\end{lemma}

\begin{vershort}
\begin{proof}
[Proof of Lemma \ref{lem.sol.perm.aj-ai-sum.3}.]If $i=j$, then all four truth
values $\left[  i<j\right]  $, $\left[  \sigma\left(  i\right)  >\sigma\left(
j\right)  \right]  $, $\left[  j<i\right]  $ and $\left[  \sigma\left(
j\right)  >\sigma\left(  i\right)  \right]  $ equal $0$ (because $i=j$ and
$\sigma\left(  i\right)  =\sigma\left(  j\right)  $), whereas the two truth
values $\left[  i\leq j\right]  $ and $\left[  \sigma\left(  i\right)
\leq\sigma\left(  j\right)  \right]  $ equal $1$ (since $i=j\leq j$ and
$\sigma\left(  i\right)  =\sigma\left(  j\right)  \leq\sigma\left(  j\right)
$). Hence, if $i=j$, then Lemma \ref{lem.sol.perm.aj-ai-sum.3} boils down to
the equality $0\cdot0-0\cdot0=1-1$, which is obvious. Thus, for the rest of
the proof of Lemma \ref{lem.sol.perm.aj-ai-sum.3}, we WLOG assume that $i\neq
j$.

The map $\sigma$ is a permutation of $\left[  n\right]  $ (since $\sigma\in
S_{n}$), and thus is bijective. Hence, $\sigma$ is injective. Thus, from
$i\neq j$, we conclude that $\sigma\left(  i\right)  \neq\sigma\left(
j\right)  $. Therefore, $\sigma\left(  i\right)  \leq\sigma\left(  j\right)  $
holds if and only if $\sigma\left(  i\right)  <\sigma\left(  j\right)  $. We
thus have the following equivalence of statements:%
\[
\left(  \sigma\left(  i\right)  \leq\sigma\left(  j\right)  \right)
\ \Longleftrightarrow\ \left(  \sigma\left(  i\right)  <\sigma\left(
j\right)  \right)  \ \Longleftrightarrow\ \left(  \sigma\left(  j\right)
>\sigma\left(  i\right)  \right)  .
\]
Therefore, $\left[  \sigma\left(  i\right)  \leq\sigma\left(  j\right)
\right]  =\left[  \sigma\left(  j\right)  >\sigma\left(  i\right)  \right]  $.

Also, $i\leq j$ holds if and only if $i<j$ (since $i\neq j$). Thus, the
statements $\left(  i\leq j\right)  $ and $\left(  i<j\right)  $ are
equivalent. Hence, $\left[  i\leq j\right]  =\left[  i<j\right]  $.

But any two statements $\mathcal{A}$ and $\mathcal{B}$ satisfy
\begin{equation}
\left[  \mathcal{A}\right]  -\left[  \mathcal{B}\right]  =\left[
\mathcal{A}\right]  \left[  \text{not }\mathcal{B}\right]  -\left[  \text{not
}\mathcal{A}\right]  \left[  \mathcal{B}\right]  .
\label{pf.lem.sol.perm.aj-ai-sum.3.short.gen}%
\end{equation}
(Indeed, if $\mathcal{A}$ and $\mathcal{B}$ are two statements, then%
\begin{align*}
&  \left[  \mathcal{A}\right]  \underbrace{\left[  \text{not }\mathcal{B}%
\right]  }_{\substack{=1-\left[  \mathcal{B}\right]  \\\text{(by Exercise
\ref{exe.iverson-prop} \textbf{(b)} (applied to }\mathcal{B}\\\text{instead of
}\mathcal{A}\text{))}}}-\underbrace{\left[  \text{not }\mathcal{A}\right]
}_{\substack{=1-\left[  \mathcal{A}\right]  \\\text{(by Exercise
\ref{exe.iverson-prop} \textbf{(b)})}}}\left[  \mathcal{B}\right] \\
&  =\left[  \mathcal{A}\right]  \left(  1-\left[  \mathcal{B}\right]  \right)
-\left(  1-\left[  \mathcal{A}\right]  \right)  \left[  \mathcal{B}\right]
=\left[  \mathcal{A}\right]  -\left[  \mathcal{A}\right]  \left[
\mathcal{B}\right]  -\left[  \mathcal{B}\right]  +\left[  \mathcal{A}\right]
\left[  \mathcal{B}\right]  =\left[  \mathcal{A}\right]  -\left[
\mathcal{B}\right]  .
\end{align*}
)

Applying (\ref{pf.lem.sol.perm.aj-ai-sum.3.short.gen}) to $\mathcal{A}=\left(
i\leq j\right)  $ and $\mathcal{B}=\left(  \sigma\left(  i\right)  \leq
\sigma\left(  j\right)  \right)  $, we obtain%
\begin{align*}
&  \left[  i\leq j\right]  -\left[  \sigma\left(  i\right)  \leq\sigma\left(
j\right)  \right] \\
&  =\underbrace{\left[  i\leq j\right]  }_{=\left[  i<j\right]  }%
\underbrace{\left[  \text{not }\sigma\left(  i\right)  \leq\sigma\left(
j\right)  \right]  }_{\substack{=\left[  \sigma\left(  i\right)
>\sigma\left(  j\right)  \right]  \\\text{(since }\left(  \text{not }%
\sigma\left(  i\right)  \leq\sigma\left(  j\right)  \right)  \text{
is}\\\text{equivalent to }\left(  \sigma\left(  i\right)  >\sigma\left(
j\right)  \right)  \text{)}}}-\underbrace{\left[  \text{not }i\leq j\right]
}_{\substack{=\left[  j<i\right]  \\\text{(since }\left(  \text{not }i\leq
j\right)  \text{ is}\\\text{equivalent to }\left(  j<i\right)  \text{)}%
}}\underbrace{\left[  \sigma\left(  i\right)  \leq\sigma\left(  j\right)
\right]  }_{=\left[  \sigma\left(  j\right)  >\sigma\left(  i\right)  \right]
}\\
&  =\left[  i<j\right]  \left[  \sigma\left(  i\right)  >\sigma\left(
j\right)  \right]  -\left[  j<i\right]  \left[  \sigma\left(  j\right)
>\sigma\left(  i\right)  \right]  .
\end{align*}
This proves Lemma \ref{lem.sol.perm.aj-ai-sum.3}.
\end{proof}
\end{vershort}

\begin{verlong}
\begin{proof}
[Proof of Lemma \ref{lem.sol.perm.aj-ai-sum.3}.]Recall that $S_{n}$ is the set
of all permutations of the set $\left\{  1,2,\ldots,n\right\}  $ (by the
definition of $S_{n}$). In other words, $S_{n}$ is the set of all permutations
of the set $\left[  n\right]  $ (since $\left[  n\right]  =\left\{
1,2,\ldots,n\right\}  $).

But $\sigma\in S_{n}$. In other words, $\sigma$ is a permutation of the set
$\left[  n\right]  $ (since $S_{n}$ is the set of all permutations of the set
$\left[  n\right]  $). In other words, $\sigma$ is a bijection $\left[
n\right]  \rightarrow\left[  n\right]  $. Thus, the map $\sigma$ is bijective,
hence injective.

We are in one of the following three cases:

\textit{Case 1:} We have $i<j$.

\textit{Case 2:} We have $i=j$.

\textit{Case 3:} We have $i>j$.

Let us first consider Case 1. In this case, we have $i<j$. Hence, $j>i$. Thus,
we don't have $j<i$. Therefore, $\left[  j<i\right]  =0$. Also, from $i<j$, we
obtain $i\leq j$, and therefore $\left[  i\leq j\right]  =1$. Also, $\left[
i<j\right]  =1$ (since $i<j$).

Also, $\left(  \sigma\left(  i\right)  >\sigma\left(  j\right)  \right)  $ and
$\left(  \text{not }\sigma\left(  i\right)  \leq\sigma\left(  j\right)
\right)  $ are two equivalent logical statements. Hence, Exercise
\ref{exe.iverson-prop} \textbf{(a)} (applied to $\mathcal{A}=\left(
\sigma\left(  i\right)  >\sigma\left(  j\right)  \right)  $ and $\mathcal{B}%
=\left(  \text{not }\sigma\left(  i\right)  \leq\sigma\left(  j\right)
\right)  $) yields $\left[  \sigma\left(  i\right)  >\sigma\left(  j\right)
\right]  =\left[  \text{not }\sigma\left(  i\right)  \leq\sigma\left(
j\right)  \right]  =1-\left[  \sigma\left(  i\right)  \leq\sigma\left(
j\right)  \right]  $ (by Exercise \ref{exe.iverson-prop} \textbf{(b)} (applied
to $\mathcal{A}=\left(  \sigma\left(  i\right)  \leq\sigma\left(  j\right)
\right)  $)).

Comparing%
\begin{align*}
&  \underbrace{\left[  i<j\right]  }_{=1}\left[  \sigma\left(  i\right)
>\sigma\left(  j\right)  \right]  -\underbrace{\left[  j<i\right]  }%
_{=0}\left[  \sigma\left(  j\right)  >\sigma\left(  i\right)  \right] \\
&  =\left[  \sigma\left(  i\right)  >\sigma\left(  j\right)  \right]
-\underbrace{0\left[  \sigma\left(  j\right)  >\sigma\left(  i\right)
\right]  }_{=0}=\left[  \sigma\left(  i\right)  >\sigma\left(  j\right)
\right]  =1-\left[  \sigma\left(  i\right)  \leq\sigma\left(  j\right)
\right]
\end{align*}
with%
\[
\underbrace{\left[  i\leq j\right]  }_{=1}-\left[  \sigma\left(  i\right)
\leq\sigma\left(  j\right)  \right]  =1-\left[  \sigma\left(  i\right)
\leq\sigma\left(  j\right)  \right]  ,
\]
we obtain $\left[  i<j\right]  \left[  \sigma\left(  i\right)  >\sigma\left(
j\right)  \right]  -\left[  j<i\right]  \left[  \sigma\left(  j\right)
>\sigma\left(  i\right)  \right]  =\left[  i\leq j\right]  -\left[
\sigma\left(  i\right)  \leq\sigma\left(  j\right)  \right]  $. Thus, Lemma
\ref{lem.sol.perm.aj-ai-sum.3} is proven in Case 1.

Let us next consider Case 2. In this case, we have $i=j$. Hence, we don't have
$i<j$. Thus, $\left[  i<j\right]  =0$. Also, $j=i$; thus, we don't have $j<i$.
Hence, $\left[  j<i\right]  =0$. Finally, from $i=j$, we obtain $i\leq j$, so
that $\left[  i\leq j\right]  =1$. Also, from $i=j$, we obtain $\sigma\left(
i\right)  =\sigma\left(  j\right)  \leq\sigma\left(  j\right)  $, so that
$\left[  \sigma\left(  i\right)  \leq\sigma\left(  j\right)  \right]  =1$.
Now, comparing%
\begin{align*}
&  \underbrace{\left[  i<j\right]  }_{=0}\left[  \sigma\left(  i\right)
>\sigma\left(  j\right)  \right]  -\underbrace{\left[  j<i\right]  }%
_{=0}\left[  \sigma\left(  j\right)  >\sigma\left(  i\right)  \right] \\
&  =0\left[  \sigma\left(  i\right)  >\sigma\left(  j\right)  \right]
-0\left[  \sigma\left(  j\right)  >\sigma\left(  i\right)  \right]  =0
\end{align*}
with%
\[
\underbrace{\left[  i\leq j\right]  }_{=1}-\underbrace{\left[  \sigma\left(
i\right)  \leq\sigma\left(  j\right)  \right]  }_{=1}=1-1=0,
\]
we obtain $\left[  i<j\right]  \left[  \sigma\left(  i\right)  >\sigma\left(
j\right)  \right]  -\left[  j<i\right]  \left[  \sigma\left(  j\right)
>\sigma\left(  i\right)  \right]  =\left[  i\leq j\right]  -\left[
\sigma\left(  i\right)  \leq\sigma\left(  j\right)  \right]  $. Thus, Lemma
\ref{lem.sol.perm.aj-ai-sum.3} is proven in Case 2.

Let us finally consider Case 3. In this case, we have $i>j$. Thus, $j<i$.
Hence, $\left[  j<i\right]  =1$. Also, we don't have $i<j$ (since $i>j$);
thus, $\left[  i<j\right]  =0$. Furthermore, we don't have $i\leq j$ (since
$i>j$); thus, $\left[  i\leq j\right]  =0$.

If we had $\sigma\left(  i\right)  =\sigma\left(  j\right)  $, then we would
have $i=j$ (since $\sigma$ is injective), which would contradict $i>j$. Thus,
we cannot have $\sigma\left(  i\right)  =\sigma\left(  j\right)  $. In other
words, the statement $\sigma\left(  i\right)  =\sigma\left(  j\right)  $ is
false. But we have the following logical equivalence:%
\begin{align*}
\left(  \sigma\left(  i\right)  \leq\sigma\left(  j\right)  \right)  \  &
\Longleftrightarrow\ \left(  \sigma\left(  i\right)  <\sigma\left(  j\right)
\text{ or }\underbrace{\sigma\left(  i\right)  =\sigma\left(  j\right)
}_{\text{This statement is false}}\right) \\
&  \Longleftrightarrow\ \left(  \sigma\left(  i\right)  <\sigma\left(
j\right)  \right)  \ \Longleftrightarrow\ \left(  \sigma\left(  j\right)
>\sigma\left(  i\right)  \right)  .
\end{align*}
Thus, $\left(  \sigma\left(  i\right)  \leq\sigma\left(  j\right)  \right)  $
and $\left(  \sigma\left(  j\right)  >\sigma\left(  i\right)  \right)  $ are
two equivalent logical statements. Hence, Exercise \ref{exe.iverson-prop}
\textbf{(a)} (applied to $\mathcal{A}=\left(  \sigma\left(  i\right)
\leq\sigma\left(  j\right)  \right)  $ and $\mathcal{B}=\left(  \sigma\left(
j\right)  >\sigma\left(  i\right)  \right)  $) yields $\left[  \sigma\left(
i\right)  \leq\sigma\left(  j\right)  \right]  =\left[  \sigma\left(
j\right)  >\sigma\left(  i\right)  \right]  $.

Now, comparing%
\begin{align*}
&  \underbrace{\left[  i<j\right]  }_{=0}\left[  \sigma\left(  i\right)
>\sigma\left(  j\right)  \right]  -\underbrace{\left[  j<i\right]  }%
_{=1}\left[  \sigma\left(  j\right)  >\sigma\left(  i\right)  \right] \\
&  =\underbrace{0\left[  \sigma\left(  i\right)  >\sigma\left(  j\right)
\right]  }_{=0}-\left[  \sigma\left(  j\right)  >\sigma\left(  i\right)
\right]  =-\left[  \sigma\left(  j\right)  >\sigma\left(  i\right)  \right]
\end{align*}
with%
\[
\underbrace{\left[  i\leq j\right]  }_{=0}-\underbrace{\left[  \sigma\left(
i\right)  \leq\sigma\left(  j\right)  \right]  }_{=\left[  \sigma\left(
j\right)  >\sigma\left(  i\right)  \right]  }=0-\left[  \sigma\left(
j\right)  >\sigma\left(  i\right)  \right]  =-\left[  \sigma\left(  j\right)
>\sigma\left(  i\right)  \right]  ,
\]
we obtain $\left[  i<j\right]  \left[  \sigma\left(  i\right)  >\sigma\left(
j\right)  \right]  -\left[  j<i\right]  \left[  \sigma\left(  j\right)
>\sigma\left(  i\right)  \right]  =\left[  i\leq j\right]  -\left[
\sigma\left(  i\right)  \leq\sigma\left(  j\right)  \right]  $. Thus, Lemma
\ref{lem.sol.perm.aj-ai-sum.3} is proven in Case 3.

We have now proven Lemma \ref{lem.sol.perm.aj-ai-sum.3} in each of the three
Cases 1, 2 and 3. Since these three Cases cover all possibilities, we thus
conclude that Lemma \ref{lem.sol.perm.aj-ai-sum.3} always holds.
\end{proof}
\end{verlong}

Solving Exercise \ref{exe.perm.aj-ai-sum} is now a mere matter of computation:

\begin{proof}
[Solution to Exercise \ref{exe.perm.aj-ai-sum}.]Let $\left[  n\right]
=\left\{  1,2,\ldots,n\right\}  $. Lemma \ref{lem.sol.perm.aj-ai-sum.2} yields%
\begin{align}
&  \sum_{\substack{1\leq i<j\leq n;\\\sigma\left(  i\right)  >\sigma\left(
j\right)  }}\left(  a_{j}-a_{i}\right) \nonumber\\
&  =\sum_{i\in\left[  n\right]  }\sum_{j\in\left[  n\right]  }%
\underbrace{\left[  i<j\right]  \left[  \sigma\left(  i\right)  >\sigma\left(
j\right)  \right]  \left(  a_{j}-a_{i}\right)  }_{=\left[  i<j\right]  \left[
\sigma\left(  i\right)  >\sigma\left(  j\right)  \right]  a_{j}-\left[
i<j\right]  \left[  \sigma\left(  i\right)  >\sigma\left(  j\right)  \right]
a_{i}}\nonumber\\
&  =\sum_{i\in\left[  n\right]  }\sum_{j\in\left[  n\right]  }\left(  \left[
i<j\right]  \left[  \sigma\left(  i\right)  >\sigma\left(  j\right)  \right]
a_{j}-\left[  i<j\right]  \left[  \sigma\left(  i\right)  >\sigma\left(
j\right)  \right]  a_{i}\right) \nonumber\\
&  =\sum_{i\in\left[  n\right]  }\sum_{j\in\left[  n\right]  }\left[
i<j\right]  \left[  \sigma\left(  i\right)  >\sigma\left(  j\right)  \right]
a_{j}-\sum_{i\in\left[  n\right]  }\sum_{j\in\left[  n\right]  }\left[
i<j\right]  \left[  \sigma\left(  i\right)  >\sigma\left(  j\right)  \right]
a_{i}. \label{sol.perm.aj-ai-sum.1}%
\end{align}

\begin{vershort}
But%
\begin{align*}
&  \underbrace{\sum_{i\in\left[  n\right]  }\sum_{j\in\left[  n\right]  }%
}_{=\sum_{j\in\left[  n\right]  }\sum_{i\in\left[  n\right]  }}\left[
i<j\right]  \left[  \sigma\left(  i\right)  >\sigma\left(  j\right)  \right]
a_{i}\\
&  =\sum_{j\in\left[  n\right]  }\sum_{i\in\left[  n\right]  }\left[
i<j\right]  \left[  \sigma\left(  i\right)  >\sigma\left(  j\right)  \right]
a_{i}=\sum_{i\in\left[  n\right]  }\sum_{j\in\left[  n\right]  }\left[
j<i\right]  \left[  \sigma\left(  j\right)  >\sigma\left(  i\right)  \right]
a_{j}%
\end{align*}
(here, we have renamed the summation indices $j$ and $i$ as $i$ and $j$, respectively).
\end{vershort}

\begin{verlong}
But
\begin{align*}
&  \underbrace{\sum_{i\in\left[  n\right]  }\sum_{j\in\left[  n\right]  }%
}_{=\sum_{j\in\left[  n\right]  }\sum_{i\in\left[  n\right]  }}\left[
i<j\right]  \left[  \sigma\left(  i\right)  >\sigma\left(  j\right)  \right]
a_{i}\\
&  =\sum_{j\in\left[  n\right]  }\sum_{i\in\left[  n\right]  }\left[
i<j\right]  \left[  \sigma\left(  i\right)  >\sigma\left(  j\right)  \right]
a_{i}\\
&  =\sum_{j\in\left[  n\right]  }\sum_{u\in\left[  n\right]  }\left[
u<j\right]  \left[  \sigma\left(  u\right)  >\sigma\left(  j\right)  \right]
a_{u}\\
&  \ \ \ \ \ \ \ \ \ \ \left(
\begin{array}
[c]{c}%
\text{here, we have renamed the summation index }i\text{ as }u\\
\text{in the inner sum}%
\end{array}
\right) \\
&  =\sum_{i\in\left[  n\right]  }\sum_{u\in\left[  n\right]  }\left[
u<i\right]  \left[  \sigma\left(  u\right)  >\sigma\left(  i\right)  \right]
a_{u}\\
&  \ \ \ \ \ \ \ \ \ \ \left(
\begin{array}
[c]{c}%
\text{here, we have renamed the summation index }j\text{ as }i\\
\text{in the outer sum}%
\end{array}
\right) \\
&  =\sum_{i\in\left[  n\right]  }\sum_{j\in\left[  n\right]  }\left[
j<i\right]  \left[  \sigma\left(  j\right)  >\sigma\left(  i\right)  \right]
a_{j}\\
&  \ \ \ \ \ \ \ \ \ \ \left(
\begin{array}
[c]{c}%
\text{here, we have renamed the summation index }u\text{ as }j\\
\text{in the inner sum}%
\end{array}
\right)  .
\end{align*}

\end{verlong}

Hence, (\ref{sol.perm.aj-ai-sum.1}) becomes%
\begin{align*}
&  \sum_{\substack{1\leq i<j\leq n;\\\sigma\left(  i\right)  >\sigma\left(
j\right)  }}\left(  a_{j}-a_{i}\right) \\
&  =\sum_{i\in\left[  n\right]  }\sum_{j\in\left[  n\right]  }\left[
i<j\right]  \left[  \sigma\left(  i\right)  >\sigma\left(  j\right)  \right]
a_{j}-\underbrace{\sum_{i\in\left[  n\right]  }\sum_{j\in\left[  n\right]
}\left[  i<j\right]  \left[  \sigma\left(  i\right)  >\sigma\left(  j\right)
\right]  a_{i}}_{=\sum_{i\in\left[  n\right]  }\sum_{j\in\left[  n\right]
}\left[  j<i\right]  \left[  \sigma\left(  j\right)  >\sigma\left(  i\right)
\right]  a_{j}}\\
&  =\sum_{i\in\left[  n\right]  }\sum_{j\in\left[  n\right]  }\left[
i<j\right]  \left[  \sigma\left(  i\right)  >\sigma\left(  j\right)  \right]
a_{j}-\sum_{i\in\left[  n\right]  }\sum_{j\in\left[  n\right]  }\left[
j<i\right]  \left[  \sigma\left(  j\right)  >\sigma\left(  i\right)  \right]
a_{j}\\
&  =\underbrace{\sum_{i\in\left[  n\right]  }\sum_{j\in\left[  n\right]  }%
}_{=\sum_{j\in\left[  n\right]  }\sum_{i\in\left[  n\right]  }}%
\underbrace{\left(  \left[  i<j\right]  \left[  \sigma\left(  i\right)
>\sigma\left(  j\right)  \right]  a_{j}-\left[  j<i\right]  \left[
\sigma\left(  j\right)  >\sigma\left(  i\right)  \right]  a_{j}\right)
}_{=\left(  \left[  i<j\right]  \left[  \sigma\left(  i\right)  >\sigma\left(
j\right)  \right]  -\left[  j<i\right]  \left[  \sigma\left(  j\right)
>\sigma\left(  i\right)  \right]  \right)  a_{j}}\\
&  =\sum_{j\in\left[  n\right]  }\sum_{i\in\left[  n\right]  }%
\underbrace{\left(  \left[  i<j\right]  \left[  \sigma\left(  i\right)
>\sigma\left(  j\right)  \right]  -\left[  j<i\right]  \left[  \sigma\left(
j\right)  >\sigma\left(  i\right)  \right]  \right)  }_{\substack{=\left[
i\leq j\right]  -\left[  \sigma\left(  i\right)  \leq\sigma\left(  j\right)
\right]  \\\text{(by Lemma \ref{lem.sol.perm.aj-ai-sum.3})}}}a_{j}\\
&  =\sum_{j\in\left[  n\right]  }\sum_{i\in\left[  n\right]  }\left(  \left[
i\leq j\right]  -\left[  \sigma\left(  i\right)  \leq\sigma\left(  j\right)
\right]  \right)  a_{j}=\sum_{j\in\left[  n\right]  }\underbrace{\left(
\sum_{i\in\left[  n\right]  }\left(  \left[  i\leq j\right]  -\left[
\sigma\left(  i\right)  \leq\sigma\left(  j\right)  \right]  \right)  \right)
}_{=\sum_{i\in\left[  n\right]  }\left[  i\leq j\right]  -\sum_{i\in\left[
n\right]  }\left[  \sigma\left(  i\right)  \leq\sigma\left(  j\right)
\right]  }a_{j}\\
&  =\sum_{j\in\left[  n\right]  }\left(  \underbrace{\sum_{i\in\left[
n\right]  }\left[  i\leq j\right]  }_{\substack{=j\\\text{(by Lemma
\ref{lem.sol.perm.aj-ai-sum.1} \textbf{(a)})}}}-\underbrace{\sum_{i\in\left[
n\right]  }\left[  \sigma\left(  i\right)  \leq\sigma\left(  j\right)
\right]  }_{\substack{=\sigma\left(  j\right)  \\\text{(by Lemma
\ref{lem.sol.perm.aj-ai-sum.1} \textbf{(b)})}}}\right)  a_{j}\\
&  =\underbrace{\sum_{j\in\left[  n\right]  }}_{=\sum_{j=1}^{n}}%
\underbrace{\left(  j-\sigma\left(  j\right)  \right)  a_{j}}_{=a_{j}\left(
j-\sigma\left(  j\right)  \right)  }=\sum_{j=1}^{n}a_{j}\left(  j-\sigma
\left(  j\right)  \right)  =\sum_{i=1}^{n}a_{i}\left(  i-\sigma\left(
i\right)  \right)
\end{align*}
(here, we have renamed the summation index $j$ as $i$). This solves Exercise
\ref{exe.perm.aj-ai-sum}.
\end{proof}

\subsection{Solution to Exercise \ref{exe.perm.pij-pii-sum}}

\begin{proof}
[Solution to Exercise \ref{exe.perm.pij-pii-sum}.]We have $\pi\in S_{n}$. In
other words, $\pi$ is a permutation of $\left\{  1,2,\ldots,n\right\}  $
(since $S_{n}$ is the set of all permutations of $\left\{  1,2,\ldots
,n\right\}  $). In other words, $\pi$ is a bijection $\left\{  1,2,\ldots
,n\right\}  \rightarrow\left\{  1,2,\ldots,n\right\}  $. Hence, we can
substitute $\pi\left(  i\right)  $ for $i$ in the sum $\sum_{i\in\left\{
1,2,\ldots,n\right\}  }i^{2}$. We thus obtain%
\begin{equation}
\sum_{i\in\left\{  1,2,\ldots,n\right\}  }i^{2}=\sum_{i\in\left\{
1,2,\ldots,n\right\}  }\left(  \pi\left(  i\right)  \right)  ^{2}.
\label{sol.perm.pij-pii-sum.1}%
\end{equation}

\textbf{(a)} Exercise \ref{exe.perm.aj-ai-sum} (applied to $\sigma=\pi$ and
$a_{k}=\pi\left(  k\right)  +k$) yields%
\begin{align*}
\sum_{\substack{1\leq i<j\leq n;\\\pi\left(  i\right)  >\pi\left(  j\right)
}}\left(  \left(  \pi\left(  j\right)  +j\right)  -\left(  \pi\left(
i\right)  +i\right)  \right)   &  =\underbrace{\sum_{i=1}^{n}}_{=\sum
_{i\in\left\{  1,2,\ldots,n\right\}  }}\underbrace{\left(  \pi\left(
i\right)  +i\right)  \left(  i-\pi\left(  i\right)  \right)  }%
_{\substack{=\left(  i+\pi\left(  i\right)  \right)  \left(  i-\pi\left(
i\right)  \right)  \\=i^{2}-\left(  \pi\left(  i\right)  \right)
^{2}\\\text{(since }\left(  x+y\right)  \left(  x-y\right)  =x^{2}%
-y^{2}\\\text{for any two numbers }x\text{ and }y\text{)}}}\\
&  =\sum_{i\in\left\{  1,2,\ldots,n\right\}  }\left(  i^{2}-\left(  \pi\left(
i\right)  \right)  ^{2}\right) \\
&  =\underbrace{\sum_{i\in\left\{  1,2,\ldots,n\right\}  }i^{2}}%
_{\substack{=\sum_{i\in\left\{  1,2,\ldots,n\right\}  }\left(  \pi\left(
i\right)  \right)  ^{2}\\\text{(by (\ref{sol.perm.pij-pii-sum.1}))}}%
}-\sum_{i\in\left\{  1,2,\ldots,n\right\}  }\left(  \pi\left(  i\right)
\right)  ^{2}\\
&  =\sum_{i\in\left\{  1,2,\ldots,n\right\}  }\left(  \pi\left(  i\right)
\right)  ^{2}-\sum_{i\in\left\{  1,2,\ldots,n\right\}  }\left(  \pi\left(
i\right)  \right)  ^{2}=0.
\end{align*}
Hence,%
\begin{align*}
0  &  =\sum_{\substack{1\leq i<j\leq n;\\\pi\left(  i\right)  >\pi\left(
j\right)  }}\underbrace{\left(  \left(  \pi\left(  j\right)  +j\right)
-\left(  \pi\left(  i\right)  +i\right)  \right)  }_{=\left(  \pi\left(
j\right)  -\pi\left(  i\right)  \right)  -\left(  i-j\right)  }\\
&  =\sum_{\substack{1\leq i<j\leq n;\\\pi\left(  i\right)  >\pi\left(
j\right)  }}\left(  \left(  \pi\left(  j\right)  -\pi\left(  i\right)
\right)  -\left(  i-j\right)  \right) \\
&  =\sum_{\substack{1\leq i<j\leq n;\\\pi\left(  i\right)  >\pi\left(
j\right)  }}\left(  \pi\left(  j\right)  -\pi\left(  i\right)  \right)
-\sum_{\substack{1\leq i<j\leq n;\\\pi\left(  i\right)  >\pi\left(  j\right)
}}\left(  i-j\right)  .
\end{align*}
Adding $\sum_{\substack{1\leq i<j\leq n;\\\pi\left(  i\right)  >\pi\left(
j\right)  }}\left(  i-j\right)  $ to both sides of this equality, we obtain%
\[
\sum_{\substack{1\leq i<j\leq n;\\\pi\left(  i\right)  >\pi\left(  j\right)
}}\left(  i-j\right)  =\sum_{\substack{1\leq i<j\leq n;\\\pi\left(  i\right)
>\pi\left(  j\right)  }}\left(  \pi\left(  j\right)  -\pi\left(  i\right)
\right)  .
\]
This solves Exercise \ref{exe.perm.pij-pii-sum} \textbf{(a)}.

\textbf{(b)} Let $w_{0}$ denote the permutation in $S_{n}$ which sends each
$k\in\left\{  1,2,\ldots,n\right\}  $ to $n+1-k$. Define a permutation
$\sigma\in S_{n}$ by $\sigma=w_{0}\circ\pi$. Thus, each $k\in\left\{
1,2,\ldots,n\right\}  $ satisfies%
\begin{equation}
\underbrace{\sigma}_{=w_{0}\circ\pi}\left(  k\right)  =\left(  w_{0}\circ
\pi\right)  \left(  k\right)  =w_{0}\left(  \pi\left(  k\right)  \right)
=n+1-\pi\left(  k\right)  \label{sol.perm.pij-pii-sum.b.1}%
\end{equation}
(by the definition of $w_{0}$).

For any $\left(  i,j\right)  \in\left\{  1,2,\ldots,n\right\}  ^{2}$, we have
the following chain of logical equivalences:%
\begin{align*}
\left(  \underbrace{\sigma\left(  i\right)  }_{\substack{=n+1-\pi\left(
i\right)  \\\text{(by (\ref{sol.perm.pij-pii-sum.b.1})}\\\text{(applied to
}k=i\text{))}}}>\underbrace{\sigma\left(  j\right)  }_{\substack{=n+1-\pi
\left(  j\right)  \\\text{(by (\ref{sol.perm.pij-pii-sum.b.1})}%
\\\text{(applied to }k=j\text{))}}}\right)  \  &  \Longleftrightarrow\ \left(
n+1-\pi\left(  i\right)  >n+1-\pi\left(  j\right)  \right) \\
&  \Longleftrightarrow\ \left(  \pi\left(  i\right)  <\pi\left(  j\right)
\right)  .
\end{align*}
Thus, for any $\left(  i,j\right)  \in\left\{  1,2,\ldots,n\right\}  ^{2}$,
the condition $\left(  \sigma\left(  i\right)  >\sigma\left(  j\right)
\right)  $ is equivalent to $\left(  \pi\left(  i\right)  <\pi\left(
j\right)  \right)  $. Hence, the summation sign \textquotedblleft%
$\sum_{\substack{1\leq i<j\leq n;\\\sigma\left(  i\right)  >\sigma\left(
j\right)  }}$\textquotedblright\ can be rewritten as \textquotedblleft%
$\sum_{\substack{1\leq i<j\leq n;\\\pi\left(  i\right)  <\pi\left(  j\right)
}}$\textquotedblright. In other words, we have%
\[
\sum_{\substack{1\leq i<j\leq n;\\\sigma\left(  i\right)  >\sigma\left(
j\right)  }}=\sum_{\substack{1\leq i<j\leq n;\\\pi\left(  i\right)
<\pi\left(  j\right)  }}
\]
(an equality between summation signs). Now, Exercise
\ref{exe.perm.pij-pii-sum} \textbf{(a)} (applied to $\sigma$ instead of $\pi$)
yields%
\begin{align*}
\sum_{\substack{1\leq i<j\leq n;\\\sigma\left(  i\right)  >\sigma\left(
j\right)  }}\left(  \sigma\left(  j\right)  -\sigma\left(  i\right)  \right)
&  =\underbrace{\sum_{\substack{1\leq i<j\leq n;\\\sigma\left(  i\right)
>\sigma\left(  j\right)  }}}_{=\sum_{\substack{1\leq i<j\leq n;\\\pi\left(
i\right)  <\pi\left(  j\right)  }}}\underbrace{\left(  i-j\right)
}_{=-\left(  j-i\right)  }=\sum_{\substack{1\leq i<j\leq n;\\\pi\left(
i\right)  <\pi\left(  j\right)  }}\left(  -\left(  j-i\right)  \right) \\
&  =-\sum_{\substack{1\leq i<j\leq n;\\\pi\left(  i\right)  <\pi\left(
j\right)  }}\left(  j-i\right)  .
\end{align*}
Comparing this with%
\begin{align*}
&  \underbrace{\sum_{\substack{1\leq i<j\leq n;\\\sigma\left(  i\right)
>\sigma\left(  j\right)  }}}_{=\sum_{\substack{1\leq i<j\leq n;\\\pi\left(
i\right)  <\pi\left(  j\right)  }}}\left(  \underbrace{\sigma\left(  j\right)
}_{\substack{=n+1-\pi\left(  j\right)  \\\text{(by
(\ref{sol.perm.pij-pii-sum.b.1})}\\\text{(applied to }k=j\text{))}%
}}-\underbrace{\sigma\left(  i\right)  }_{\substack{=n+1-\pi\left(  i\right)
\\\text{(by (\ref{sol.perm.pij-pii-sum.b.1})}\\\text{(applied to }%
k=i\text{))}}}\right) \\
&  =\sum_{\substack{1\leq i<j\leq n;\\\pi\left(  i\right)  <\pi\left(
j\right)  }}\underbrace{\left(  \left(  n+1-\pi\left(  j\right)  \right)
-\left(  n+1-\pi\left(  i\right)  \right)  \right)  }_{=-\left(  \pi\left(
j\right)  -\pi\left(  i\right)  \right)  }=\sum_{\substack{1\leq i<j\leq
n;\\\pi\left(  i\right)  <\pi\left(  j\right)  }}\left(  -\left(  \pi\left(
j\right)  -\pi\left(  i\right)  \right)  \right) \\
&  =-\sum_{\substack{1\leq i<j\leq n;\\\pi\left(  i\right)  <\pi\left(
j\right)  }}\left(  \pi\left(  j\right)  -\pi\left(  i\right)  \right)  ,
\end{align*}
we obtain%
\[
-\sum_{\substack{1\leq i<j\leq n;\\\pi\left(  i\right)  <\pi\left(  j\right)
}}\left(  \pi\left(  j\right)  -\pi\left(  i\right)  \right)  =-\sum
_{\substack{1\leq i<j\leq n;\\\pi\left(  i\right)  <\pi\left(  j\right)
}}\left(  j-i\right)  .
\]
Thus,%
\[
\sum_{\substack{1\leq i<j\leq n;\\\pi\left(  i\right)  <\pi\left(  j\right)
}}\left(  \pi\left(  j\right)  -\pi\left(  i\right)  \right)  =\sum
_{\substack{1\leq i<j\leq n;\\\pi\left(  i\right)  <\pi\left(  j\right)
}}\left(  j-i\right)  .
\]
This solves Exercise \ref{exe.perm.pij-pii-sum} \textbf{(b)}.
\end{proof}

\subsection{\label{sect.sols.perm.rearrangement}Solution to Exercise
\ref{exe.perm.rearrangement}}

Throughout Section \ref{sect.sols.perm.rearrangement}, we shall use the
notations introduced in Definition \ref{def.set12...m} and in Definition
\ref{def.transpos.ii}.

Before we step to the solution of Exercise \ref{exe.perm.rearrangement}, let
us prove a simple lemma:

\begin{lemma}
\label{lem.sol.perm.rearrangement.tpi}Let $n\in\mathbb{N}$. Let $\pi\in S_{n}%
$. Let $u$ and $v$ be two elements of $\left[  n\right]  $. Then,%
\[
t_{\pi\left(  u\right)  ,\pi\left(  v\right)  }\circ\pi=\pi\circ t_{u,v}.
\]

\end{lemma}

\begin{proof}
[Proof of Lemma \ref{lem.sol.perm.rearrangement.tpi}.]We have $\left[
n\right]  =\left\{  1,2,\ldots,n\right\}  $ (by the definition of $\left[
n\right]  $).

\begin{vershort}
Thus, $S_{n}$ is the set of all permutations of $\left[  n\right]  $. Hence,
all three maps $t_{\pi\left(  u\right)  ,\pi\left(  v\right)  }$, $\pi$ and
$t_{u,v}$ are permutations of $\left[  n\right]  $, that is, bijections from
$\left[  n\right]  $ to $\left[  n\right]  $. Thus, the map $\pi$ is a
bijection, and hence is injective.
\end{vershort}

\begin{verlong}
Hence, $u$ and $v$ are elements of $\left\{  1,2,\ldots,n\right\}  $ (since
$u$ and $v$ are elements of $\left[  n\right]  $). Thus, the permutation
$t_{u,v}$ in $S_{n}$ is well-defined.

Recall that $S_{n}$ is the set of all permutations of $\left\{  1,2,\ldots
,n\right\}  $. In other words, $S_{n}$ is the set of all permutations of
$\left[  n\right]  $ (since $\left[  n\right]  =\left\{  1,2,\ldots,n\right\}
$). But $\pi$ belongs to $S_{n}$. Thus, $\pi$ is a permutation of $\left[
n\right]  $ (since $S_{n}$ is the set of all permutations of $\left[
n\right]  $). In other words, $\pi$ is a bijection from $\left[  n\right]  $
to $\left[  n\right]  $. Hence, $\pi\left(  u\right)  $ and $\pi\left(
v\right)  $ are two elements of $\left[  n\right]  $. In other words,
$\pi\left(  u\right)  $ and $\pi\left(  v\right)  $ are two elements of
$\left\{  1,2,\ldots,n\right\}  $ (since $\left[  n\right]  =\left\{
1,2,\ldots,n\right\}  $). Thus, the permutation $t_{\pi\left(  u\right)
,\pi\left(  v\right)  }$ in $S_{n}$ is well-defined.

Both $t_{\pi\left(  u\right)  ,\pi\left(  v\right)  }$ and $t_{u,v}$ belong to
$S_{n}$, and thus are permutations of $\left[  n\right]  $ (since $S_{n}$ is
the set of all permutations of $\left[  n\right]  $). In other words, both
$t_{\pi\left(  u\right)  ,\pi\left(  v\right)  }$ and $t_{u,v}$ are bijections
from $\left[  n\right]  $ to $\left[  n\right]  $.

The map $\pi$ is a bijection, i.e., a bijective map. Thus, $\pi$ is injective
and surjective.

Recall that $\pi$, $t_{\pi\left(  u\right)  ,\pi\left(  v\right)  }$ and
$t_{u,v}$ are bijections from $\left[  n\right]  $ to $\left[  n\right]  $,
therefore maps from $\left[  n\right]  $ to $\left[  n\right]  $. Hence, both
$t_{\pi\left(  u\right)  ,\pi\left(  v\right)  }\circ\pi$ and $\pi\circ
t_{u,v}$ are maps from $\left[  n\right]  $ to $\left[  n\right]  $.
\end{verlong}

Fix $k\in\left[  n\right]  $. We shall prove that $\left(  t_{\pi\left(
u\right)  ,\pi\left(  v\right)  }\circ\pi\right)  \left(  k\right)  =\left(
\pi\circ t_{u,v}\right)  \left(  k\right)  $.

We are in one of the following three cases:

\textit{Case 1:} We have $k=u$.

\textit{Case 2:} We have $k=v$.

\textit{Case 3:} We have neither $k=u$ nor $k=v$.

Let us first consider Case 1. In this case, we have $k=u$. But Lemma
\ref{lem.sol.transpose.code.tij1} \textbf{(a)} (applied to $i=u$ and $j=v$)
yields $t_{u,v}\left(  u\right)  =v$. Also,%
\[
\left(  t_{\pi\left(  u\right)  ,\pi\left(  v\right)  }\circ\pi\right)
\left(  \underbrace{k}_{=u}\right)  =\left(  t_{\pi\left(  u\right)
,\pi\left(  v\right)  }\circ\pi\right)  \left(  u\right)  =t_{\pi\left(
u\right)  ,\pi\left(  v\right)  }\left(  \pi\left(  u\right)  \right)
=\pi\left(  v\right)
\]
(by Lemma \ref{lem.sol.transpose.code.tij1} \textbf{(a)} (applied to
$i=\pi\left(  u\right)  $ and $j=\pi\left(  v\right)  $)). Comparing this with%
\[
\left(  \pi\circ t_{u,v}\right)  \left(  \underbrace{k}_{=u}\right)  =\left(
\pi\circ t_{u,v}\right)  \left(  u\right)  =\pi\left(  \underbrace{t_{u,v}%
\left(  u\right)  }_{=v}\right)  =\pi\left(  v\right)  ,
\]
we obtain $\left(  t_{\pi\left(  u\right)  ,\pi\left(  v\right)  }\circ
\pi\right)  \left(  k\right)  =\left(  \pi\circ t_{u,v}\right)  \left(
k\right)  $. Hence, $\left(  t_{\pi\left(  u\right)  ,\pi\left(  v\right)
}\circ\pi\right)  \left(  k\right)  =\left(  \pi\circ t_{u,v}\right)  \left(
k\right)  $ is proven in Case 1.

\begin{vershort}
The argument in Case 2 is analogous, and we leave it to the reader.
\end{vershort}

\begin{verlong}
Let us next consider Case 2. In this case, we have $k=v$. But Lemma
\ref{lem.sol.transpose.code.tij1} \textbf{(b)} (applied to $i=u$ and $j=v$)
yields $t_{u,v}\left(  v\right)  =u$. Also,%
\[
\left(  t_{\pi\left(  u\right)  ,\pi\left(  v\right)  }\circ\pi\right)
\left(  \underbrace{k}_{=v}\right)  =\left(  t_{\pi\left(  u\right)
,\pi\left(  v\right)  }\circ\pi\right)  \left(  v\right)  =t_{\pi\left(
u\right)  ,\pi\left(  v\right)  }\left(  \pi\left(  v\right)  \right)
=\pi\left(  u\right)
\]
(by Lemma \ref{lem.sol.transpose.code.tij1} \textbf{(b)} (applied to
$i=\pi\left(  u\right)  $ and $j=\pi\left(  v\right)  $)). Comparing this with%
\[
\left(  \pi\circ t_{u,v}\right)  \left(  \underbrace{k}_{=v}\right)  =\left(
\pi\circ t_{u,v}\right)  \left(  v\right)  =\pi\left(  \underbrace{t_{u,v}%
\left(  v\right)  }_{=u}\right)  =\pi\left(  u\right)  ,
\]
we obtain $\left(  t_{\pi\left(  u\right)  ,\pi\left(  v\right)  }\circ
\pi\right)  \left(  k\right)  =\left(  \pi\circ t_{u,v}\right)  \left(
k\right)  $. Hence, $\left(  t_{\pi\left(  u\right)  ,\pi\left(  v\right)
}\circ\pi\right)  \left(  k\right)  =\left(  \pi\circ t_{u,v}\right)  \left(
k\right)  $ is proven in Case 2.
\end{verlong}

\begin{vershort}
Let us now consider Case 3. In this case, we have neither $k=u$ nor $k=v$.
Thus, $k\in\left[  n\right]  \setminus\left\{  u,v\right\}  $. Hence, Lemma
\ref{lem.sol.transpose.code.tij1} \textbf{(c)} (applied to $i=u$ and $j=v$)
yields $t_{u,v}\left(  k\right)  =k$. Also, recall that we have neither $k=u$
nor $k=v$. Thus, we have neither $\pi\left(  k\right)  =\pi\left(  u\right)  $
nor $\pi\left(  k\right)  =\pi\left(  v\right)  $ (since the map $\pi$ is
injective). In other words, $\pi\left(  k\right)  \in\left[  n\right]
\setminus\left\{  \pi\left(  u\right)  ,\pi\left(  v\right)  \right\}  $.
Thus, Lemma \ref{lem.sol.transpose.code.tij1} \textbf{(c)} (applied to
$\pi\left(  u\right)  $, $\pi\left(  v\right)  $ and $\pi\left(  k\right)  $
instead of $i$, $j$ and $k$) yields $t_{\pi\left(  u\right)  ,\pi\left(
v\right)  }\left(  \pi\left(  k\right)  \right)  =\pi\left(  k\right)  $. Now,%
\[
\left(  t_{\pi\left(  u\right)  ,\pi\left(  v\right)  }\circ\pi\right)
\left(  k\right)  =t_{\pi\left(  u\right)  ,\pi\left(  v\right)  }\left(
\pi\left(  k\right)  \right)  =\pi\left(  k\right)  .
\]
Comparing this with%
\[
\left(  \pi\circ t_{u,v}\right)  \left(  k\right)  =\pi\left(
\underbrace{t_{u,v}\left(  k\right)  }_{=k}\right)  =\pi\left(  k\right)  ,
\]
we obtain $\left(  t_{\pi\left(  u\right)  ,\pi\left(  v\right)  }\circ
\pi\right)  \left(  k\right)  =\left(  \pi\circ t_{u,v}\right)  \left(
k\right)  $. Hence, $\left(  t_{\pi\left(  u\right)  ,\pi\left(  v\right)
}\circ\pi\right)  \left(  k\right)  =\left(  \pi\circ t_{u,v}\right)  \left(
k\right)  $ is proven in Case 3.
\end{vershort}

\begin{verlong}
Let us now consider Case 3. In this case, we have neither $k=u$ nor $k=v$.
Thus, $k\notin\left\{  u,v\right\}  $\ \ \ \ \footnote{\textit{Proof.} Assume
the contrary. Thus, $k\in\left\{  u,v\right\}  $. Hence, either $k=u$ or
$k=v$. This contradicts the fact that we have neither $k=u$ nor $k=v$. This
contradiction shows that our assumption was false. Qed.}. Combining
$k\in\left[  n\right]  $ with $k\notin\left\{  u,v\right\}  $, we obtain
$k\in\left[  n\right]  \setminus\left\{  u,v\right\}  $. Hence, Lemma
\ref{lem.sol.transpose.code.tij1} \textbf{(c)} (applied to $i=u$ and $j=v$)
yields $t_{u,v}\left(  k\right)  =k$.

Also, recall that we have neither $k=u$ nor $k=v$. Thus, we have neither
$\pi\left(  k\right)  =\pi\left(  u\right)  $ nor $\pi\left(  k\right)
=\pi\left(  v\right)  $\ \ \ \ \footnote{\textit{Proof.} We have neither $k=u$
nor $k=v$. Hence, $k\neq u$ and $k\neq v$.
\par
If we had $\pi\left(  k\right)  =\pi\left(  u\right)  $, then we would have
$k=u$ (since the map $\pi$ is injective), which would contradict $k\neq u$.
Thus, we cannot have $\pi\left(  k\right)  =\pi\left(  u\right)  $. In other
words, we have $\pi\left(  k\right)  \neq\pi\left(  u\right)  $.
\par
If we had $\pi\left(  k\right)  =\pi\left(  v\right)  $, then we would have
$k=v$ (since the map $\pi$ is injective), which would contradict $k\neq v$.
Thus, we cannot have $\pi\left(  k\right)  =\pi\left(  v\right)  $. In other
words, we have $\pi\left(  k\right)  \neq\pi\left(  v\right)  $.
\par
Thus, we have $\pi\left(  k\right)  \neq\pi\left(  u\right)  $ and $\pi\left(
k\right)  \neq\pi\left(  v\right)  $. In other words, we have neither
$\pi\left(  k\right)  =\pi\left(  u\right)  $ nor $\pi\left(  k\right)
=\pi\left(  v\right)  $.}. Thus, $\pi\left(  k\right)  \notin\left\{
\pi\left(  u\right)  ,\pi\left(  v\right)  \right\}  $%
\ \ \ \ \footnote{\textit{Proof.} Assume the contrary. Thus, $\pi\left(
k\right)  \in\left\{  \pi\left(  u\right)  ,\pi\left(  v\right)  \right\}  $.
Hence, either $\pi\left(  k\right)  =\pi\left(  u\right)  $ or $\pi\left(
k\right)  =\pi\left(  v\right)  $. This contradicts the fact that we have
neither $\pi\left(  k\right)  =\pi\left(  u\right)  $ nor $\pi\left(
k\right)  =\pi\left(  v\right)  $. This contradiction shows that our
assumption was false. Qed.}. But $\pi\left(  k\right)  \in\left[  n\right]  $
(since $\pi$ is a bijection from $\left[  n\right]  $ to $\left[  n\right]
$). Combining this with $\pi\left(  k\right)  \notin\left\{  \pi\left(
u\right)  ,\pi\left(  v\right)  \right\}  $, we obtain $\pi\left(  k\right)
\in\left[  n\right]  \setminus\left\{  \pi\left(  u\right)  ,\pi\left(
v\right)  \right\}  $. Hence, Lemma \ref{lem.sol.transpose.code.tij1}
\textbf{(c)} (applied to $\pi\left(  u\right)  $, $\pi\left(  v\right)  $ and
$\pi\left(  k\right)  $ instead of $i$, $j$ and $k$) yields $t_{\pi\left(
u\right)  ,\pi\left(  v\right)  }\left(  \pi\left(  k\right)  \right)
=\pi\left(  k\right)  $. Now,%
\[
\left(  t_{\pi\left(  u\right)  ,\pi\left(  v\right)  }\circ\pi\right)
\left(  k\right)  =t_{\pi\left(  u\right)  ,\pi\left(  v\right)  }\left(
\pi\left(  k\right)  \right)  =\pi\left(  k\right)  .
\]
Comparing this with%
\[
\left(  \pi\circ t_{u,v}\right)  \left(  k\right)  =\pi\left(
\underbrace{t_{u,v}\left(  k\right)  }_{=k}\right)  =\pi\left(  k\right)  ,
\]
we obtain $\left(  t_{\pi\left(  u\right)  ,\pi\left(  v\right)  }\circ
\pi\right)  \left(  k\right)  =\left(  \pi\circ t_{u,v}\right)  \left(
k\right)  $. Hence, $\left(  t_{\pi\left(  u\right)  ,\pi\left(  v\right)
}\circ\pi\right)  \left(  k\right)  =\left(  \pi\circ t_{u,v}\right)  \left(
k\right)  $ is proven in Case 3.
\end{verlong}

We have now proven $\left(  t_{\pi\left(  u\right)  ,\pi\left(  v\right)
}\circ\pi\right)  \left(  k\right)  =\left(  \pi\circ t_{u,v}\right)  \left(
k\right)  $ in each of the three Cases 1, 2 and 3. Thus, $\left(
t_{\pi\left(  u\right)  ,\pi\left(  v\right)  }\circ\pi\right)  \left(
k\right)  =\left(  \pi\circ t_{u,v}\right)  \left(  k\right)  $ always holds.

Forget now that we fixed $k$. We thus have shown that $\left(  t_{\pi\left(
u\right)  ,\pi\left(  v\right)  }\circ\pi\right)  \left(  k\right)  =\left(
\pi\circ t_{u,v}\right)  \left(  k\right)  $ for each $k\in\left[  n\right]
$. In other words, $t_{\pi\left(  u\right)  ,\pi\left(  v\right)  }\circ
\pi=\pi\circ t_{u,v}$ (since both $t_{\pi\left(  u\right)  ,\pi\left(
v\right)  }\circ\pi$ and $\pi\circ t_{u,v}$ are maps from $\left[  n\right]  $
to $\left[  n\right]  $). Thus, Lemma \ref{lem.sol.perm.rearrangement.tpi} is proven.
\end{proof}

We can now solve Exercise \ref{exe.perm.rearrangement}:

\begin{proof}
[Solution to Exercise \ref{exe.perm.rearrangement}.]Recall that $\left(
i_{1},i_{2},\ldots,i_{n}\right)  \in\left[  1\right]  \times\left[  2\right]
\times\cdots\times\left[  n\right]  $. Thus,%
\begin{equation}
i_{j}\in\left[  j\right]  \ \ \ \ \ \ \ \ \ \ \text{for each }j\in\left[
n\right]  . \label{sol.perm.rearrangement.1}%
\end{equation}

\begin{verlong}
Hence,%
\begin{equation}
i_{j}\in\left[  n\right]  \ \ \ \ \ \ \ \ \ \ \text{for each }j\in\left[
n\right]  \label{sol.perm.rearrangement.1b}%
\end{equation}
\footnote{\textit{Proof of (\ref{sol.perm.rearrangement.1b}):} Let
$j\in\left[  n\right]  $. Then, $j\in\left[  n\right]  =\left\{
1,2,\ldots,n\right\}  $, so that $j\leq n$ and thus $\left\{  1,2,\ldots
,j\right\}  \subseteq\left\{  1,2,\ldots,n\right\}  $. Now,
(\ref{sol.perm.rearrangement.1}) yields $i_{j}\in\left[  j\right]  =\left\{
1,2,\ldots,j\right\}  \subseteq\left\{  1,2,\ldots,n\right\}  =\left[
n\right]  $. This proves (\ref{sol.perm.rearrangement.1b}).}.
\end{verlong}

The definition of $\sigma_{0}$ shows that%
\begin{equation}
\sigma_{0}=t_{1,i_{1}}\circ t_{2,i_{2}}\circ\cdots\circ t_{0,i_{0}}=\left(
\text{a composition of }0\text{ permutations}\right)  =\operatorname*{id}.
\label{sol.perm.rearrangement.sigma0}%
\end{equation}
The definition of $\sigma_{n}$ shows that%
\begin{equation}
\sigma_{n}=t_{1,i_{1}}\circ t_{2,i_{2}}\circ\cdots\circ t_{n,i_{n}}%
=\sigma\label{sol.perm.rearrangement.sigman}%
\end{equation}
(since $\sigma=t_{1,i_{1}}\circ t_{2,i_{2}}\circ\cdots\circ t_{n,i_{n}}$).

\begin{verlong}
We have $\left[  n\right]  =\left\{  1,2,\ldots,n\right\}  $ (by the
definition of $\left[  n\right]  $). Recall that $S_{n}$ is the set of all
permutations of $\left\{  1,2,\ldots,n\right\}  $. In other words, $S_{n}$ is
the set of all permutations of $\left[  n\right]  $ (since $\left[  n\right]
=\left\{  1,2,\ldots,n\right\}  $). But $\sigma$ belongs to $S_{n}$. Thus,
$\sigma$ is a permutation of $\left[  n\right]  $ (since $S_{n}$ is the set of
all permutations of $\left[  n\right]  $). In other words, $\sigma$ is a
bijection from $\left[  n\right]  $ to $\left[  n\right]  $. Moreover, for
each $k\in\left\{  0,1,\ldots,n\right\}  $, the map $\sigma_{k}$ belongs to
$S_{n}$, and thus is a permutation of $\left[  n\right]  $ (since $S_{n}$ is
the set of all permutations of $\left[  n\right]  $). In other words, for each
$k\in\left\{  0,1,\ldots,n\right\}  $, the map $\sigma_{k}$ is a bijection
from $\left[  n\right]  $ to $\left[  n\right]  $.
\end{verlong}

For each $k\in\left[  n\right]  $, we have%
\begin{equation}
\sigma_{k}=\sigma_{k-1}\circ t_{k,i_{k}}. \label{sol.perm.rearrangement.rec1}%
\end{equation}

\begin{vershort}
[\textit{Proof of (\ref{sol.perm.rearrangement.rec1}):} Let $k\in\left[
n\right]  $. The definition of $\sigma_{k-1}$ yields $\sigma_{k-1}=t_{1,i_{1}%
}\circ t_{2,i_{2}}\circ\cdots\circ t_{k-1,i_{k-1}}$. But the definition of
$\sigma_{k}$ yields%
\[
\sigma_{k}=t_{1,i_{1}}\circ t_{2,i_{2}}\circ\cdots\circ t_{k,i_{k}%
}=\underbrace{\left(  t_{1,i_{1}}\circ t_{2,i_{2}}\circ\cdots\circ
t_{k-1,i_{k-1}}\right)  }_{=\sigma_{k-1}}\circ t_{k,i_{k}}=\sigma_{k-1}\circ
t_{k,i_{k}}.
\]
This proves (\ref{sol.perm.rearrangement.rec1}).]
\end{vershort}

\begin{verlong}
[\textit{Proof of (\ref{sol.perm.rearrangement.rec1}):} Let $k\in\left[
n\right]  $. Thus, $k\in\left[  n\right]  =\left\{  1,2,\ldots,n\right\}  $,
so that $k\geq1$. Also, from $k\in\left\{  1,2,\ldots,n\right\}  $, we obtain
$k-1\in\left\{  0,1,\ldots,n-1\right\}  \subseteq\left\{  0,1,\ldots
,n\right\}  $. Hence, $\sigma_{k-1}$ is well-defined. The definition of
$\sigma_{k-1}$ yields%
\begin{equation}
\sigma_{k-1}=t_{1,i_{1}}\circ t_{2,i_{2}}\circ\cdots\circ t_{k-1,i_{k-1}}.
\label{sol.perm.rearrangement.rec1.pf.1}%
\end{equation}
But $k\in\left\{  1,2,\ldots,n\right\}  \subseteq\left\{  0,1,\ldots
,n\right\}  $. Hence, $\sigma_{k}$ is well-defined as well. The definition of
$\sigma_{k}$ yields%
\begin{align*}
\sigma_{k}  &  =t_{1,i_{1}}\circ t_{2,i_{2}}\circ\cdots\circ t_{k,i_{k}%
}=\underbrace{\left(  t_{1,i_{1}}\circ t_{2,i_{2}}\circ\cdots\circ
t_{k-1,i_{k-1}}\right)  }_{=\sigma_{k-1}}\circ t_{k,i_{k}}%
\ \ \ \ \ \ \ \ \ \ \left(  \text{since }k\geq1\right) \\
&  =\sigma_{k-1}\circ t_{k,i_{k}}.
\end{align*}
This proves (\ref{sol.perm.rearrangement.rec1}).]
\end{verlong}

\textbf{(a)} Let $i\in\left[  n\right]  $. We shall prove that%
\begin{equation}
\sigma_{k}\left(  i\right)  =i\ \ \ \ \ \ \ \ \ \ \text{for each }k\in\left\{
0,1,\ldots,i-1\right\}  . \label{sol.perm.rearrangement.a.1}%
\end{equation}

[\textit{Proof of (\ref{sol.perm.rearrangement.a.1}):} We shall prove
(\ref{sol.perm.rearrangement.a.1}) by induction on $k$:

\textit{Induction base:} From (\ref{sol.perm.rearrangement.sigma0}), we obtain
$\sigma_{0}=\operatorname*{id}$. Hence, $\sigma_{0}\left(  i\right)
=\operatorname*{id}\left(  i\right)  =i$. In other words,
(\ref{sol.perm.rearrangement.a.1}) holds for $k=0$. This completes the
induction base.

\textit{Induction step:} Let $p\in\left\{  1,2,\ldots,i-1\right\}  $. Assume
that (\ref{sol.perm.rearrangement.a.1}) holds for $k=p-1$. We must prove that
(\ref{sol.perm.rearrangement.a.1}) holds for $k=p$.

\begin{verlong}
We have $p\in\left\{  1,2,\ldots,i-1\right\}  $, thus $p-1\in\left\{
0,1,\ldots,\left(  i-1\right)  -1\right\}  \subseteq\left\{  0,1,\ldots
,i-1\right\}  $. Also, $p\in\left\{  1,2,\ldots,i-1\right\}  \subseteq\left\{
1,2,\ldots,n\right\}  $ (since $i-1\leq i\leq n$ (because $i\in\left[
n\right]  =\left\{  1,2,\ldots,n\right\}  $)). Thus, $p\in\left\{
1,2,\ldots,n\right\}  =\left[  n\right]  $. Hence,
(\ref{sol.perm.rearrangement.1b}) (applied to $j=p$) yields $i_{p}\in\left[
n\right]  $.
\end{verlong}

We assumed that (\ref{sol.perm.rearrangement.a.1}) holds for $k=p-1$. In other
words, $\sigma_{p-1}\left(  i\right)  =i$. Also,
(\ref{sol.perm.rearrangement.rec1}) (applied to $k=p$) yields $\sigma
_{p}=\sigma_{p-1}\circ t_{p,i_{p}}$ (since $p\in\left[  n\right]  $).

We have $p\in\left\{  1,2,\ldots,i-1\right\}  $, thus $p\leq i-1<i$ and
therefore $p\neq i$. In other words, $i\neq p$. Also,
(\ref{sol.perm.rearrangement.1}) (applied to $j=p$) yields $i_{p}\in\left[
p\right]  =\left\{  1,2,\ldots,p\right\}  $; thus, $i_{p}\leq p<i$ and
therefore $i_{p}\neq i$. In other words, $i\neq i_{p}$. So we know that $i\neq
p$ and $i\neq i_{p}$.

\begin{vershort}
Hence, $i\in\left[  n\right]  \setminus\left\{  p,i_{p}\right\}  $. Thus,
Lemma \ref{lem.sol.transpose.code.tij1} \textbf{(c)} (applied to $p$, $i_{p}$
and $i$ instead of $i$, $j$ and $k$) yields $t_{p,i_{p}}\left(  i\right)  =i$.
\end{vershort}

\begin{verlong}
In other words, we have neither $i=p$ nor $i=i_{p}$. Hence, $i\notin\left\{
p,i_{p}\right\}  $\ \ \ \ \footnote{\textit{Proof.} Assume the contrary. Thus,
$i\in\left\{  p,i_{p}\right\}  $. In other words, we have either $i=p$ or
$i=i_{p}$. This contradicts the fact that we have neither $i=p$ nor $i=i_{p}$.
This contradiction shows that our assumption was false. Qed.}. Combining
$i\in\left[  n\right]  $ with $i\notin\left\{  p,i_{p}\right\}  $, we obtain
$i\in\left[  n\right]  \setminus\left\{  p,i_{p}\right\}  $. Thus, Lemma
\ref{lem.sol.transpose.code.tij1} \textbf{(c)} (applied to $p$, $i_{p}$ and
$i$ instead of $i$, $j$ and $k$) yields $t_{p,i_{p}}\left(  i\right)  =i$
(since $p\in\left[  n\right]  $ and $i_{p}\in\left[  n\right]  $).
\end{verlong}

Now,%
\[
\underbrace{\sigma_{p}}_{=\sigma_{p-1}\circ t_{p,i_{p}}}\left(  i\right)
=\left(  \sigma_{p-1}\circ t_{p,i_{p}}\right)  \left(  i\right)  =\sigma
_{p-1}\left(  \underbrace{t_{p,i_{p}}\left(  i\right)  }_{=i}\right)
=\sigma_{p-1}\left(  i\right)  =i.
\]
In other words, (\ref{sol.perm.rearrangement.a.1}) holds for $k=p$. This
completes the induction step. Thus, (\ref{sol.perm.rearrangement.a.1}) is
proven by induction.]

Thus, we have shown that $\sigma_{k}\left(  i\right)  =i$ for each
$k\in\left\{  0,1,\ldots,i-1\right\}  $. This solves Exercise
\ref{exe.perm.rearrangement} \textbf{(a)}.

\textbf{(b)} Let $k\in\left[  n\right]  $. The definition of $m_{k}$ yields
$m_{k}=\sigma_{k}\left(  k\right)  $.

\begin{verlong}
We have $k\in\left[  n\right]  =\left\{  1,2,\ldots,n\right\}  \subseteq
\mathbb{N}$. Recall that $\sigma_{k}\in S_{n}$; in other words, $\sigma_{k}$
is a permutation of $\left[  n\right]  $ (since $S_{n}$ is the set of all
permutations of $\left[  n\right]  $). In other words, $\sigma_{k}$ is a
bijection from $\left[  n\right]  $ to $\left[  n\right]  $. Hence, the map
$\sigma_{k}$ is bijective, and therefore injective. Also, $\sigma_{k}\left(
k\right)  \in\left[  n\right]  $ (since $\sigma_{k}$ is a bijection from
$\left[  n\right]  $ to $\left[  n\right]  $).
\end{verlong}

Assume (for the sake of contradiction) that $m_{k}\notin\left[  k\right]  $.
Combining $m_{k}=\sigma_{k}\left(  k\right)  \in\left[  n\right]  $ with
$m_{k}\notin\left[  k\right]  $, we obtain
\[
m_{k}\in\underbrace{\left[  n\right]  }_{=\left\{  1,2,\ldots,n\right\}
}\setminus\underbrace{\left[  k\right]  }_{=\left\{  1,2,\ldots,k\right\}
}=\left\{  1,2,\ldots,n\right\}  \setminus\left\{  1,2,\ldots,k\right\}
=\left\{  k+1,k+2,\ldots,n\right\}  .
\]
Thus, $m_{k}\geq k+1$, so that $k\leq m_{k}-1$. Hence, $k\in\left\{
0,1,\ldots,m_{k}-1\right\}  $ (because $k\in\mathbb{N}$). Hence, Exercise
\ref{exe.perm.rearrangement} \textbf{(a)} (applied to $i=m_{k}$) yields
$\sigma_{k}\left(  m_{k}\right)  =m_{k}=\sigma_{k}\left(  k\right)  $. Since
$\sigma_{k}$ is injective (because $\sigma_{k}\in S_{n}$), this yields
$m_{k}=k\leq m_{k}-1<m_{k}$. But this is clearly absurd. This contradiction
shows that our assumption (that $m_{k}\notin\left[  k\right]  $) was false.
Hence, we must have $m_{k}\in\left[  k\right]  $. This solves Exercise
\ref{exe.perm.rearrangement} \textbf{(b)}.

\textbf{(c)} Let $k\in\left[  n\right]  $.

\begin{vershort}
Thus, $k\geq1$, and therefore $k-1\in\mathbb{N}$, so that $k-1\in\left\{
0,1,\ldots,k-1\right\}  $.
\end{vershort}

\begin{verlong}
We have $k\in\left[  n\right]  =\left\{  1,2,\ldots,n\right\}  $, so that
$k-1\in\left\{  0,1,\ldots,n-1\right\}  \subseteq\mathbb{N}$. Hence,
$k-1\in\left\{  0,1,\ldots,k-1\right\}  $.
\end{verlong}

Hence, Exercise \ref{exe.perm.rearrangement} \textbf{(a)} (applied to $k$ and
$k-1$ instead of $i$ and $k$) yields
\begin{equation}
\sigma_{k-1}\left(  k\right)  =k. \label{sol.perm.rearrangement.d.0}%
\end{equation}

\begin{vershort}
On the other hand, Lemma \ref{lem.sol.transpose.code.tij1} \textbf{(b)}
(applied to $k$ and $i_{k}$ instead of $i$ and $j$) yields $t_{k,i_{k}}\left(
i_{k}\right)  =k$. Furthermore, Lemma \ref{lem.sol.transpose.code.tij1}
\textbf{(a)} (applied to $k$ and $i_{k}$ instead of $i$ and $j$) yields
$t_{k,i_{k}}\left(  k\right)  =i_{k}$.
\end{vershort}

\begin{verlong}
On the other hand, (\ref{sol.perm.rearrangement.1b}) (applied to $j=k$) yields
$i_{k}\in\left[  n\right]  $. Hence, Lemma \ref{lem.sol.transpose.code.tij1}
\textbf{(b)} (applied to $k$ and $i_{k}$ instead of $i$ and $j$) yields
$t_{k,i_{k}}\left(  i_{k}\right)  =k$. Furthermore, Lemma
\ref{lem.sol.transpose.code.tij1} \textbf{(a)} (applied to $k$ and $i_{k}$
instead of $i$ and $j$) yields $t_{k,i_{k}}\left(  k\right)  =i_{k}$.
\end{verlong}

But (\ref{sol.perm.rearrangement.rec1}) yields $\sigma_{k}=\sigma_{k-1}\circ
t_{k,i_{k}}$. Hence,%
\[
\underbrace{\sigma_{k}}_{=\sigma_{k-1}\circ t_{k,i_{k}}}\left(  i_{k}\right)
=\left(  \sigma_{k-1}\circ t_{k,i_{k}}\right)  \left(  i_{k}\right)
=\sigma_{k-1}\left(  \underbrace{t_{k,i_{k}}\left(  i_{k}\right)  }%
_{=k}\right)  =\sigma_{k-1}\left(  k\right)  =k.
\]
This solves Exercise \ref{exe.perm.rearrangement} \textbf{(c)}.

\textbf{(d)} Let $k\in\left[  n\right]  $. Then,
(\ref{sol.perm.rearrangement.rec1}) yields $\sigma_{k}=\sigma_{k-1}\circ
t_{k,i_{k}}$. In our above solution to Exercise \ref{exe.perm.rearrangement}
\textbf{(c)}, we have shown that $\sigma_{k-1}\left(  k\right)  =k$ and
$t_{k,i_{k}}\left(  k\right)  =i_{k}$. Also, the definition of $m_{k}$ yields%
\[
m_{k}=\underbrace{\sigma_{k}}_{=\sigma_{k-1}\circ t_{k,i_{k}}}\left(
k\right)  =\left(  \sigma_{k-1}\circ t_{k,i_{k}}\right)  \left(  k\right)
=\sigma_{k-1}\left(  \underbrace{t_{k,i_{k}}\left(  k\right)  }_{=i_{k}%
}\right)  =\sigma_{k-1}\left(  i_{k}\right)  ,
\]
so that
\begin{equation}
\sigma_{k-1}\left(  i_{k}\right)  =m_{k}. \label{sol.perm.rearrangement.d.1}%
\end{equation}

Define $\pi\in S_{n}$ by $\pi=\sigma_{k-1}$. Then, $\underbrace{\pi}%
_{=\sigma_{k-1}}\left(  k\right)  =\sigma_{k-1}\left(  k\right)  =k$ and
$\underbrace{\pi}_{=\sigma_{k-1}}\left(  i_{k}\right)  =\sigma_{k-1}\left(
i_{k}\right)  =m_{k}$.

\begin{verlong}
On the other hand, (\ref{sol.perm.rearrangement.1b}) (applied to $j=k$) yields
$i_{k}\in\left[  n\right]  $.
\end{verlong}

Now, Lemma \ref{lem.sol.perm.rearrangement.tpi} (applied to $u=k$ and
$v=i_{k}$) yields
\[
t_{\pi\left(  k\right)  ,\pi\left(  i_{k}\right)  }\circ\pi=\underbrace{\pi
}_{=\sigma_{k-1}}\circ t_{k,i_{k}}=\sigma_{k-1}\circ t_{k,i_{k}}=\sigma_{k}%
\]
(since $\sigma_{k}=\sigma_{k-1}\circ t_{k,i_{k}}$). In view of $\pi\left(
k\right)  =k$ and $\pi\left(  i_{k}\right)  =m_{k}$, this rewrites as
$t_{k,m_{k}}\circ\pi=\sigma_{k}$. Hence, $\sigma_{k}=t_{k,m_{k}}%
\circ\underbrace{\pi}_{=\sigma_{k-1}}=t_{k,m_{k}}\circ\sigma_{k-1}$. This
solves Exercise \ref{exe.perm.rearrangement} \textbf{(d)}.

\textbf{(e)} We shall show that%
\begin{equation}
\sigma_{p}^{-1}=t_{1,m_{1}}\circ t_{2,m_{2}}\circ\cdots\circ t_{p,m_{p}%
}\ \ \ \ \ \ \ \ \ \ \text{for each }p\in\left\{  0,1,\ldots,n\right\}  .
\label{sol.perm.rearrangement.e.goal}%
\end{equation}

[\textit{Proof of (\ref{sol.perm.rearrangement.e.goal}):} We shall prove
(\ref{sol.perm.rearrangement.e.goal}) by induction on $p$:

\textit{Induction base:} We have $\sigma_{0}=\operatorname*{id}$ (by
(\ref{sol.perm.rearrangement.sigma0})) and thus $\sigma_{0}^{-1}%
=\operatorname*{id}\nolimits^{-1}=\operatorname*{id}$. Comparing this with%
\[
t_{1,m_{1}}\circ t_{2,m_{2}}\circ\cdots\circ t_{0,m_{0}}=\left(  \text{a
composition of }0\text{ permutations}\right)  =\operatorname*{id},
\]
we obtain $\sigma_{0}^{-1}=t_{1,m_{1}}\circ t_{2,m_{2}}\circ\cdots\circ
t_{0,m_{0}}$. In other words, (\ref{sol.perm.rearrangement.e.goal}) holds for
$p=0$. This completes the induction base.

\textit{Induction step:} Let $k\in\left\{  1,2,\ldots,n\right\}  $. Assume
that (\ref{sol.perm.rearrangement.e.goal}) holds for $p=k-1$. We must prove
that (\ref{sol.perm.rearrangement.e.goal}) holds for $p=k$.

We have assumed that (\ref{sol.perm.rearrangement.e.goal}) holds for $p=k-1$.
That is, we have%
\[
\sigma_{k-1}^{-1}=t_{1,m_{1}}\circ t_{2,m_{2}}\circ\cdots\circ t_{k-1,m_{k-1}%
}.
\]

\begin{vershort}
Lemma \ref{lem.sol.transpose.code.tij1} \textbf{(d)} (applied to $k$ and
$m_{k}$ instead of $i$ and $j$) yields $t_{k,m_{k}}\circ t_{k,m_{k}%
}=\operatorname*{id}$. Thus, $t_{k,m_{k}}^{-1}=t_{k,m_{k}}$.
\end{vershort}

\begin{verlong}
The definition of $m_{k}$ yields $m_{k}=\sigma_{k}\left(  k\right)  \in\left[
n\right]  $. Hence, Lemma \ref{lem.sol.transpose.code.tij1} \textbf{(d)}
(applied to $k$ and $m_{k}$ instead of $i$ and $j$) yields $t_{k,m_{k}}\circ
t_{k,m_{k}}=\operatorname*{id}$. Thus, $t_{k,m_{k}}^{-1}=t_{k,m_{k}}$.
\end{verlong}

But $k\in\left\{  1,2,\ldots,n\right\}  =\left[  n\right]  $. Hence, Exercise
\ref{exe.perm.rearrangement} \textbf{(d)} yields $\sigma_{k}=t_{k,m_{k}}%
\circ\sigma_{k-1}$. Thus,%
\begin{align*}
\sigma_{k}^{-1}  &  =\left(  t_{k,m_{k}}\circ\sigma_{k-1}\right)
^{-1}=\underbrace{\sigma_{k-1}^{-1}}_{=t_{1,m_{1}}\circ t_{2,m_{2}}\circ
\cdots\circ t_{k-1,m_{k-1}}}\circ\underbrace{t_{k,m_{k}}^{-1}}_{=t_{k,m_{k}}%
}\\
&  =\left(  t_{1,m_{1}}\circ t_{2,m_{2}}\circ\cdots\circ t_{k-1,m_{k-1}%
}\right)  \circ t_{k,m_{k}}=t_{1,m_{1}}\circ t_{2,m_{2}}\circ\cdots\circ
t_{k,m_{k}}.
\end{align*}
In other words, (\ref{sol.perm.rearrangement.e.goal}) holds for $p=k$. This
completes the induction step. Thus, (\ref{sol.perm.rearrangement.e.goal}) is
proven by induction.]

Now, $n\in\left\{  0,1,\ldots,n\right\}  $. Hence,
(\ref{sol.perm.rearrangement.e.goal}) (applied to $p=n$) yields $\sigma
_{n}^{-1}=t_{1,m_{1}}\circ t_{2,m_{2}}\circ\cdots\circ t_{n,m_{n}}$. In view
of (\ref{sol.perm.rearrangement.sigman}), this rewrites as $\sigma
^{-1}=t_{1,m_{1}}\circ t_{2,m_{2}}\circ\cdots\circ t_{n,m_{n}}$. This solves
Exercise \ref{exe.perm.rearrangement} \textbf{(e)}.

\textbf{(f)} For each permutation $\tau\in S_{n}$, we define a number
$z\left(  \tau\right)  $ by
\[
z\left(  \tau\right)  =\sum_{k=1}^{n}x_{k}y_{\tau\left(  k\right)  }.
\]

We shall show that
\begin{equation}
z\left(  \sigma_{p-1}\right)  -z\left(  \sigma_{p}\right)  =\left(  x_{i_{p}%
}-x_{p}\right)  \left(  y_{m_{p}}-y_{p}\right)  \ \ \ \ \ \ \ \ \ \ \text{for
each }p\in\left[  n\right]  . \label{sol.perm.rearrangement.f.goal}%
\end{equation}

\begin{vershort}
[\textit{Proof of (\ref{sol.perm.rearrangement.f.goal}):} Let $p\in\left[
n\right]  $. Applying (\ref{sol.perm.rearrangement.rec1}) to $k=p$, we obtain
$\sigma_{p}=\sigma_{p-1}\circ t_{p,i_{p}}$. Hence, if $p=i_{p}$, then
(\ref{sol.perm.rearrangement.f.goal}) holds\footnote{\textit{Proof.} Assume
that $p=i_{p}$. Thus, $i_{p}=p$, so that $x_{i_{p}}-x_{p}=x_{p}-x_{p}=0$.
Hence, the right hand side of (\ref{sol.perm.rearrangement.f.goal}) equals
$0$. Also, $\sigma_{p}=\sigma_{p-1}\circ\underbrace{t_{p,i_{p}}}%
_{\substack{=\operatorname*{id}\\\text{(since }p=i_{p}\text{)}}}=\sigma_{p-1}%
$, so that $z\left(  \sigma_{p-1}\right)  -z\left(  \sigma_{p}\right)
=z\left(  \sigma_{p-1}\right)  -z\left(  \sigma_{p-1}\right)  =0$. Thus, the
left hand side of (\ref{sol.perm.rearrangement.f.goal}) equals $0$ as well.
Hence, the equality (\ref{sol.perm.rearrangement.f.goal}) holds (since both
its right hand side and its left hand side equal $0$).}. Thus, for the rest of
this proof of (\ref{sol.perm.rearrangement.f.goal}), we WLOG assume that
$p\neq i_{p}$. Hence, $t_{p,i_{p}}$ is an actual transposition (not the
identity map).

From $\sigma_{p}=\sigma_{p-1}\circ t_{p,i_{p}}$, we obtain%
\[
\sigma_{p}\left(  p\right)  =\left(  \sigma_{p-1}\circ t_{p,i_{p}}\right)
\left(  p\right)  =\sigma_{p-1}\left(  \underbrace{t_{p,i_{p}}\left(
p\right)  }_{=i_{p}}\right)  =\sigma_{p-1}\left(  i_{p}\right)  ,
\]
so that%
\begin{equation}
\sigma_{p-1}\left(  i_{p}\right)  =\sigma_{p}\left(  p\right)  =m_{p}
\label{sol.perm.rearrangement.f.goal.pf.short.1}%
\end{equation}
(since the definition of $m_{p}$ yields $m_{p}=\sigma_{p}\left(  p\right)  $).

From $\sigma_{p}=\sigma_{p-1}\circ t_{p,i_{p}}$, we also obtain%
\[
\sigma_{p}\left(  i_{p}\right)  =\left(  \sigma_{p-1}\circ t_{p,i_{p}}\right)
\left(  i_{p}\right)  =\sigma_{p-1}\left(  \underbrace{t_{p,i_{p}}\left(
i_{p}\right)  }_{=p}\right)  =\sigma_{p-1}\left(  p\right)  ,
\]
so that%
\begin{equation}
\sigma_{p-1}\left(  p\right)  =\sigma_{p}\left(  i_{p}\right)  =p
\label{sol.perm.rearrangement.f.goal.pf.short.2}%
\end{equation}
(by Exercise \ref{exe.perm.rearrangement} \textbf{(c)}, applied to $k=p$).

Every $k\in\left[  n\right]  $ satisfying $k\neq p$ and $k\neq i_{p}$
satisfies
\begin{equation}
\sigma_{p-1}\left(  k\right)  =\sigma_{p}\left(  k\right)
\label{sol.perm.rearrangement.f.goal.pf.short.3}%
\end{equation}
\footnote{\textit{Proof:} Let $k\in\left[  n\right]  $ be such that $k\neq p$
and $k\neq i_{p}$. Thus, $t_{p,i_{p}}\left(  k\right)  =k$. But $\sigma
_{p}=\sigma_{p-1}\circ t_{p,i_{p}}$; hence, $\sigma_{p}\left(  k\right)
=\left(  \sigma_{p-1}\circ t_{p,i_{p}}\right)  \left(  k\right)  =\sigma
_{p-1}\left(  \underbrace{t_{p,i_{p}}\left(  k\right)  }_{=k}\right)
=\sigma_{p-1}\left(  k\right)  $, so that $\sigma_{p-1}\left(  k\right)
=\sigma_{p}\left(  k\right)  $, qed.}. Now, the definition of $z\left(
\sigma_{p-1}\right)  $ yields%
\begin{align*}
z\left(  \sigma_{p-1}\right)   &  =\sum_{k=1}^{n}x_{k}y_{\sigma_{p-1}\left(
k\right)  }=x_{p}\underbrace{y_{\sigma_{p-1}\left(  p\right)  }}%
_{\substack{=y_{p}\\\text{(by (\ref{sol.perm.rearrangement.f.goal.pf.short.2}%
))}}}+x_{i_{p}}\underbrace{y_{\sigma_{p-1}\left(  i_{p}\right)  }%
}_{\substack{=y_{m_{p}}\\\text{(by
(\ref{sol.perm.rearrangement.f.goal.pf.short.1}))}}}+\sum_{\substack{k\in
\left[  n\right]  ;\\k\neq p\text{ and }k\neq i_{p}}}x_{k}%
\underbrace{y_{\sigma_{p-1}\left(  k\right)  }}_{\substack{=y_{\sigma
_{p}\left(  k\right)  }\\\text{(by
(\ref{sol.perm.rearrangement.f.goal.pf.short.3}))}}}\\
&  \ \ \ \ \ \ \ \ \ \ \left(
\begin{array}
[c]{c}%
\text{here, we have split the addends for }k=p\text{ and}\\
\text{for }k=i_{p}\text{ from the sum (and these are}\\
\text{two distinct addends, since }p\neq i_{p}\text{)}%
\end{array}
\right) \\
&  =x_{p}y_{p}+x_{i_{p}}y_{m_{p}}+\sum_{\substack{k\in\left[  n\right]
;\\k\neq p\text{ and }k\neq i_{p}}}x_{k}y_{\sigma_{p}\left(  k\right)  }.
\end{align*}
On the other hand, the definition of $z\left(  \sigma_{p}\right)  $ yields%
\begin{align*}
z\left(  \sigma_{p}\right)   &  =\sum_{k=1}^{n}x_{k}y_{\sigma_{p}\left(
k\right)  }=x_{p}\underbrace{y_{\sigma_{p}\left(  p\right)  }}%
_{\substack{=y_{m_{p}}\\\text{(since }\sigma_{p}\left(  p\right)
=m_{p}\text{)}}}+x_{i_{p}}\underbrace{y_{\sigma_{p}\left(  i_{p}\right)  }%
}_{\substack{=y_{p}\\\text{(since }\sigma_{p}\left(  i_{p}\right)  =p\text{)}%
}}+\sum_{\substack{k\in\left[  n\right]  ;\\k\neq p\text{ and }k\neq i_{p}%
}}x_{k}y_{\sigma_{p}\left(  k\right)  }\\
&  \ \ \ \ \ \ \ \ \ \ \left(
\begin{array}
[c]{c}%
\text{here, we have split the addends for }k=p\text{ and}\\
\text{for }k=i_{p}\text{ from the sum (and these are}\\
\text{two distinct addends, since }p\neq i_{p}\text{)}%
\end{array}
\right) \\
&  =x_{p}y_{m_{p}}+x_{i_{p}}y_{p}+\sum_{\substack{k\in\left[  n\right]
;\\k\neq p\text{ and }k\neq i_{p}}}x_{k}y_{\sigma_{p}\left(  k\right)  }.
\end{align*}
Subtracting this equality from the preceding equality, we obtain%
\begin{align*}
&  z\left(  \sigma_{p-1}\right)  -z\left(  \sigma_{p}\right) \\
&  =\left(  x_{p}y_{p}+x_{i_{p}}y_{m_{p}}+\sum_{\substack{k\in\left[
n\right]  ;\\k\neq p\text{ and }k\neq i_{p}}}x_{k}y_{\sigma_{p}\left(
k\right)  }\right)  -\left(  x_{p}y_{m_{p}}+x_{i_{p}}y_{p}+\sum
_{\substack{k\in\left[  n\right]  ;\\k\neq p\text{ and }k\neq i_{p}}%
}x_{k}y_{\sigma_{p}\left(  k\right)  }\right) \\
&  =x_{p}y_{p}+x_{i_{p}}y_{m_{p}}-x_{p}y_{m_{p}}-x_{i_{p}}y_{p}=\left(
x_{i_{p}}-x_{p}\right)  \left(  y_{m_{p}}-y_{p}\right)  .
\end{align*}
This proves (\ref{sol.perm.rearrangement.f.goal}).]
\end{vershort}

\begin{verlong}
[\textit{Proof of (\ref{sol.perm.rearrangement.f.goal}):} Let $p\in\left[
n\right]  $. Applying (\ref{sol.perm.rearrangement.rec1}) to $k=p$, we obtain
$\sigma_{p}=\sigma_{p-1}\circ t_{p,i_{p}}$. Hence, if $p=i_{p}$, then
(\ref{sol.perm.rearrangement.f.goal}) holds\footnote{\textit{Proof.} Assume
that $p=i_{p}$. We must prove that (\ref{sol.perm.rearrangement.f.goal})
holds.
\par
We have $p=i_{p}$. Thus, $i_{p}=p$, so that $x_{i_{p}}-x_{p}=x_{p}-x_{p}=0$.
\par
On the other hand, Definition \ref{def.transpos.ii} yields $t_{p,i_{p}%
}=\operatorname*{id}$ (since $p=i_{p}$). Hence, $\sigma_{p}=\sigma_{p-1}%
\circ\underbrace{t_{p,i_{p}}}_{=\operatorname*{id}}=\sigma_{p-1}$, so that
$z\left(  \sigma_{p-1}\right)  -z\left(  \sigma_{p}\right)  =z\left(
\sigma_{p-1}\right)  -z\left(  \sigma_{p-1}\right)  =0$. Comparing this with
$\underbrace{\left(  x_{i_{p}}-x_{p}\right)  }_{=0}\left(  y_{m_{p}}%
-y_{p}\right)  =0$, we obtain $z\left(  \sigma_{p-1}\right)  -z\left(
\sigma_{p}\right)  =\left(  x_{i_{p}}-x_{p}\right)  \left(  y_{m_{p}}%
-y_{p}\right)  $. Hence, the equality (\ref{sol.perm.rearrangement.f.goal})
holds. Qed.}. Thus, for the rest of this proof of
(\ref{sol.perm.rearrangement.f.goal}), we can WLOG assume that we don't have
$p=i_{p}$. Assume this.

We have $p\neq i_{p}$ (since we don't have $p=i_{p}$).

Now, (\ref{sol.perm.rearrangement.1b}) (applied to $j=p$) yields $i_{p}%
\in\left[  n\right]  $. Hence, $p$ and $i_{p}$ are two elements of $\left[
n\right]  $. Thus, $\left\{  p,i_{p}\right\}  \subseteq\left[  n\right]  $.

Exercise \ref{exe.perm.rearrangement} \textbf{(c)} (applied to $k=p$) yields
\[
\sigma_{p}\left(  i_{p}\right)  =p.
\]

The definition of $m_{p}$ yields $m_{p}=\sigma_{p}\left(  p\right)  $; thus,%
\[
\sigma_{p}\left(  p\right)  =m_{p}.
\]

Recall the equality (\ref{sol.perm.rearrangement.d.1}) that we have shown (for
each $k\in\left[  n\right]  $) during our solution to Exercise
\ref{exe.perm.rearrangement} \textbf{(d)}. Applying this equality to $k=p$, we
obtain%
\begin{equation}
\sigma_{p-1}\left(  i_{p}\right)  =m_{p}.
\label{sol.perm.rearrangement.f.goal.pf.1}%
\end{equation}

Also, recall the equality (\ref{sol.perm.rearrangement.d.0}) that we have
shown (for each $k\in\left[  n\right]  $) during our solution to Exercise
\ref{exe.perm.rearrangement} \textbf{(c)}. Applying this equality to $k=p$, we
obtain%
\begin{equation}
\sigma_{p-1}\left(  p\right)  =p. \label{sol.perm.rearrangement.f.goal.pf.2}%
\end{equation}

Every $k\in\left[  n\right]  \setminus\left\{  p,i_{p}\right\}  $ satisfies
\begin{equation}
\sigma_{p-1}\left(  k\right)  =\sigma_{p}\left(  k\right)
\label{sol.perm.rearrangement.f.goal.pf.3}%
\end{equation}
\footnote{\textit{Proof:} Let $k\in\left[  n\right]  \setminus\left\{
p,i_{p}\right\}  $. Thus, Lemma \ref{lem.sol.transpose.code.tij1} \textbf{(c)}
(applied to $p$ and $i_{p}$ instead of $i$ and $j$) yields $t_{p,i_{p}}\left(
k\right)  =k$. But $\sigma_{p}=\sigma_{p-1}\circ t_{p,i_{p}}$; hence,
$\sigma_{p}\left(  k\right)  =\left(  \sigma_{p-1}\circ t_{p,i_{p}}\right)
\left(  k\right)  =\sigma_{p-1}\left(  \underbrace{t_{p,i_{p}}\left(
k\right)  }_{=k}\right)  =\sigma_{p-1}\left(  k\right)  $, so that
$\sigma_{p-1}\left(  k\right)  =\sigma_{p}\left(  k\right)  $. This proves
(\ref{sol.perm.rearrangement.f.goal.pf.3}).}. Now, the definition of $z\left(
\sigma_{p-1}\right)  $ yields%
\begin{align}
z\left(  \sigma_{p-1}\right)   &  =\underbrace{\sum_{k=1}^{n}}%
_{\substack{=\sum_{k\in\left[  n\right]  }\\\text{(since }\left[  n\right]
=\left\{  1,2,\ldots,n\right\}  \text{)}}}x_{k}y_{\sigma_{p-1}\left(
k\right)  }=\sum_{k\in\left[  n\right]  }x_{k}y_{\sigma_{p-1}\left(  k\right)
}\nonumber\\
&  =\underbrace{\sum_{\substack{k\in\left[  n\right]  ;\\k\in\left\{
p,i_{p}\right\}  }}}_{\substack{=\sum_{k\in\left\{  p,i_{p}\right\}
}\\\text{(since }\left\{  p,i_{p}\right\}  \subseteq\left[  n\right]
\text{)}}}x_{k}y_{\sigma_{p-1}\left(  k\right)  }+\underbrace{\sum
_{\substack{k\in\left[  n\right]  ;\\k\notin\left\{  p,i_{p}\right\}  }%
}}_{=\sum_{k\in\left[  n\right]  \setminus\left\{  p,i_{p}\right\}  }}%
x_{k}y_{\sigma_{p-1}\left(  k\right)  }\nonumber\\
&  \ \ \ \ \ \ \ \ \ \ \left(
\begin{array}
[c]{c}%
\text{since each }k\in\left[  n\right]  \text{ satisfies either }k\in\left\{
p,i_{p}\right\}  \text{ or }k\notin\left\{  p,i_{p}\right\} \\
\text{(but not both at the same time)}%
\end{array}
\right) \nonumber\\
&  =\underbrace{\sum_{k\in\left\{  p,i_{p}\right\}  }x_{k}y_{\sigma
_{p-1}\left(  k\right)  }}_{\substack{=x_{p}y_{\sigma_{p-1}\left(  p\right)
}+x_{i_{p}}y_{\sigma_{p-1}\left(  i_{p}\right)  }\\\text{(since }p\neq
i_{p}\text{)}}}+\sum_{k\in\left[  n\right]  \setminus\left\{  p,i_{p}\right\}
}x_{k}\underbrace{y_{\sigma_{p-1}\left(  k\right)  }}_{\substack{=y_{\sigma
_{p}\left(  k\right)  }\\\text{(since }\sigma_{p-1}\left(  k\right)
=\sigma_{p}\left(  k\right)  \\\text{(by
(\ref{sol.perm.rearrangement.f.goal.pf.3})))}}}\nonumber\\
&  =x_{p}\underbrace{y_{\sigma_{p-1}\left(  p\right)  }}_{\substack{=y_{p}%
\\\text{(since }\sigma_{p-1}\left(  p\right)  =p\\\text{(by
(\ref{sol.perm.rearrangement.f.goal.pf.2})))}}}+x_{i_{p}}\underbrace{y_{\sigma
_{p-1}\left(  i_{p}\right)  }}_{\substack{=y_{m_{p}}\\\text{(since }%
\sigma_{p-1}\left(  i_{p}\right)  =m_{p}\\\text{(by
(\ref{sol.perm.rearrangement.f.goal.pf.1})))}}}+\sum_{k\in\left[  n\right]
\setminus\left\{  p,i_{p}\right\}  }x_{k}y_{\sigma_{p}\left(  k\right)
}\nonumber\\
&  =x_{p}y_{p}+x_{i_{p}}y_{m_{p}}+\sum_{k\in\left[  n\right]  \setminus
\left\{  p,i_{p}\right\}  }x_{k}y_{\sigma_{p}\left(  k\right)  }.
\label{sol.perm.rearrangement.f.goal.pf.p-1}%
\end{align}
On the other hand, the definition of $z\left(  \sigma_{p}\right)  $ yields%
\begin{align}
z\left(  \sigma_{p}\right)   &  =\underbrace{\sum_{k=1}^{n}}_{\substack{=\sum
_{k\in\left[  n\right]  }\\\text{(since }\left[  n\right]  =\left\{
1,2,\ldots,n\right\}  \text{)}}}x_{k}y_{\sigma_{p}\left(  k\right)  }%
=\sum_{k\in\left[  n\right]  }x_{k}y_{\sigma_{p}\left(  k\right)  }\nonumber\\
&  =\underbrace{\sum_{\substack{k\in\left[  n\right]  ;\\k\in\left\{
p,i_{p}\right\}  }}}_{\substack{=\sum_{k\in\left\{  p,i_{p}\right\}
}\\\text{(since }\left\{  p,i_{p}\right\}  \subseteq\left[  n\right]
\text{)}}}x_{k}y_{\sigma_{p}\left(  k\right)  }+\underbrace{\sum
_{\substack{k\in\left[  n\right]  ;\\k\notin\left\{  p,i_{p}\right\}  }%
}}_{=\sum_{k\in\left[  n\right]  \setminus\left\{  p,i_{p}\right\}  }}%
x_{k}y_{\sigma_{p}\left(  k\right)  }\nonumber\\
&  \ \ \ \ \ \ \ \ \ \ \left(
\begin{array}
[c]{c}%
\text{since each }k\in\left[  n\right]  \text{ satisfies either }k\in\left\{
p,i_{p}\right\}  \text{ or }k\notin\left\{  p,i_{p}\right\} \\
\text{(but not both at the same time)}%
\end{array}
\right) \nonumber\\
&  =\underbrace{\sum_{k\in\left\{  p,i_{p}\right\}  }x_{k}y_{\sigma_{p}\left(
k\right)  }}_{\substack{=x_{p}y_{\sigma_{p}\left(  p\right)  }+x_{i_{p}%
}y_{\sigma_{p}\left(  i_{p}\right)  }\\\text{(since }p\neq i_{p}\text{)}%
}}+\sum_{k\in\left[  n\right]  \setminus\left\{  p,i_{p}\right\}  }%
x_{k}y_{\sigma_{p}\left(  k\right)  }\nonumber\\
&  =x_{p}\underbrace{y_{\sigma_{p}\left(  p\right)  }}_{\substack{=y_{m_{p}%
}\\\text{(since }\sigma_{p}\left(  p\right)  =m_{p}\text{)}}}+x_{i_{p}%
}\underbrace{y_{\sigma_{p}\left(  i_{p}\right)  }}_{\substack{=y_{p}%
\\\text{(since }\sigma_{p}\left(  i_{p}\right)  =p\text{)}}}+\sum_{k\in\left[
n\right]  \setminus\left\{  p,i_{p}\right\}  }x_{k}y_{\sigma_{p}\left(
k\right)  }\nonumber\\
&  =x_{p}y_{m_{p}}+x_{i_{p}}y_{p}+\sum_{k\in\left[  n\right]  \setminus
\left\{  p,i_{p}\right\}  }x_{k}y_{\sigma_{p}\left(  k\right)  }.
\label{sol.perm.rearrangement.f.goal.pf.p}%
\end{align}
Subtracting this equality from (\ref{sol.perm.rearrangement.f.goal.pf.p-1}),
we obtain%
\begin{align*}
&  z\left(  \sigma_{p-1}\right)  -z\left(  \sigma_{p}\right) \\
&  =\left(  x_{p}y_{p}+x_{i_{p}}y_{m_{p}}+\sum_{k\in\left[  n\right]
\setminus\left\{  p,i_{p}\right\}  }x_{k}y_{\sigma_{p}\left(  k\right)
}\right)  -\left(  x_{p}y_{m_{p}}+x_{i_{p}}y_{p}+\sum_{k\in\left[  n\right]
\setminus\left\{  p,i_{p}\right\}  }x_{k}y_{\sigma_{p}\left(  k\right)
}\right) \\
&  =x_{p}y_{p}+x_{i_{p}}y_{m_{p}}-x_{p}y_{m_{p}}-x_{i_{p}}y_{p}=\left(
x_{i_{p}}-x_{p}\right)  \left(  y_{m_{p}}-y_{p}\right)  .
\end{align*}
This proves (\ref{sol.perm.rearrangement.f.goal}).]
\end{verlong}

Recall that $x_{1},x_{2},\ldots,x_{n},y_{1},y_{2},\ldots,y_{n}$ are numbers,
i.e., elements of $\mathbb{A}$, where $\mathbb{A}$ is either $\mathbb{Q}$ or
$\mathbb{R}$ or $\mathbb{C}$. Consider this $\mathbb{A}$. We have $1-1=0\leq
n$. Hence, (\ref{eq.sum.telescope}) (applied to $u=1$, $v=n$ and
$a_{s}=z\left(  \sigma_{s}\right)  $) yields%
\[
\sum_{s=1}^{n}\left(  z\left(  \sigma_{s}\right)  -z\left(  \sigma
_{s-1}\right)  \right)  =z\left(  \underbrace{\sigma_{n}}_{\substack{=\sigma
\\\text{(by (\ref{sol.perm.rearrangement.sigman}))}}}\right)  -z\left(
\underbrace{\sigma_{1-1}}_{\substack{=\sigma_{0}=\operatorname*{id}\\\text{(by
(\ref{sol.perm.rearrangement.sigma0}))}}}\right)  =z\left(  \sigma\right)
-z\left(  \operatorname*{id}\right)  .
\]
But%
\begin{align*}
&  \sum_{p=1}^{n}\left(  z\left(  \sigma_{p-1}\right)  -z\left(  \sigma
_{p}\right)  \right) \\
&  =\sum_{s=1}^{n}\underbrace{\left(  z\left(  \sigma_{s-1}\right)  -z\left(
\sigma_{s}\right)  \right)  }_{=-\left(  z\left(  \sigma_{s}\right)  -z\left(
\sigma_{s-1}\right)  \right)  }\ \ \ \ \ \ \ \ \ \ \left(
\begin{array}
[c]{c}%
\text{here, we have renamed the}\\
\text{summation index }p\text{ as }s\text{ in the sum}%
\end{array}
\right) \\
&  =\sum_{s=1}^{n}\left(  -\left(  z\left(  \sigma_{s}\right)  -z\left(
\sigma_{s-1}\right)  \right)  \right)  =-\underbrace{\sum_{s=1}^{n}\left(
z\left(  \sigma_{s}\right)  -z\left(  \sigma_{s-1}\right)  \right)
}_{=z\left(  \sigma\right)  -z\left(  \operatorname*{id}\right)  }\\
&  =-\left(  z\left(  \sigma\right)  -z\left(  \operatorname*{id}\right)
\right)  =\underbrace{z\left(  \operatorname*{id}\right)  }_{\substack{=\sum
_{k=1}^{n}x_{k}y_{\operatorname*{id}\left(  k\right)  }\\\text{(by the
definition of }z\left(  \operatorname*{id}\right)  \text{)}}%
}-\underbrace{z\left(  \sigma\right)  }_{\substack{=\sum_{k=1}^{n}%
x_{k}y_{\sigma\left(  k\right)  }\\\text{(by the definition of }z\left(
\sigma\right)  \text{)}}}\\
&  =\sum_{k=1}^{n}x_{k}\underbrace{y_{\operatorname*{id}\left(  k\right)  }%
}_{\substack{=y_{k}\\\text{(since }\operatorname*{id}\left(  k\right)
=k\text{)}}}-\sum_{k=1}^{n}x_{k}y_{\sigma\left(  k\right)  }=\sum_{k=1}%
^{n}x_{k}y_{k}-\sum_{k=1}^{n}x_{k}y_{\sigma\left(  k\right)  }.
\end{align*}

\begin{vershort}
Hence,%
\begin{align*}
&  \sum_{k=1}^{n}x_{k}y_{k}-\sum_{k=1}^{n}x_{k}y_{\sigma\left(  k\right)  }\\
&  =\sum_{p=1}^{n}\underbrace{\left(  z\left(  \sigma_{p-1}\right)  -z\left(
\sigma_{p}\right)  \right)  }_{\substack{=\left(  x_{i_{p}}-x_{p}\right)
\left(  y_{m_{p}}-y_{p}\right)  \\\text{(by
(\ref{sol.perm.rearrangement.f.goal}))}}}=\sum_{p=1}^{n}\left(  x_{i_{p}%
}-x_{p}\right)  \left(  y_{m_{p}}-y_{p}\right)  =\sum_{k=1}^{n}\left(
x_{i_{k}}-x_{k}\right)  \left(  y_{m_{k}}-y_{k}\right)
\end{align*}
(here, we have renamed the summation index $p$ as $k$). This solves Exercise
\ref{exe.perm.rearrangement} \textbf{(f)}.
\end{vershort}

\begin{verlong}
Hence,%
\begin{align*}
&  \sum_{k=1}^{n}x_{k}y_{k}-\sum_{k=1}^{n}x_{k}y_{\sigma\left(  k\right)  }\\
&  =\underbrace{\sum_{p=1}^{n}}_{\substack{=\sum_{p\in\left[  n\right]
}\\\text{(since }\left[  n\right]  =\left\{  1,2,\ldots,n\right\}  \text{)}%
}}\left(  z\left(  \sigma_{p-1}\right)  -z\left(  \sigma_{p}\right)  \right)
=\sum_{p\in\left[  n\right]  }\underbrace{\left(  z\left(  \sigma
_{p-1}\right)  -z\left(  \sigma_{p}\right)  \right)  }_{\substack{=\left(
x_{i_{p}}-x_{p}\right)  \left(  y_{m_{p}}-y_{p}\right)  \\\text{(by
(\ref{sol.perm.rearrangement.f.goal}))}}}\\
&  =\underbrace{\sum_{p\in\left[  n\right]  }}_{\substack{=\sum_{p=1}%
^{n}\\\text{(since }\left[  n\right]  =\left\{  1,2,\ldots,n\right\}
\text{)}}}\left(  x_{i_{p}}-x_{p}\right)  \left(  y_{m_{p}}-y_{p}\right)
=\sum_{p=1}^{n}\left(  x_{i_{p}}-x_{p}\right)  \left(  y_{m_{p}}-y_{p}\right)
\\
&  =\sum_{k=1}^{n}\left(  x_{i_{k}}-x_{k}\right)  \left(  y_{m_{k}}%
-y_{k}\right)
\end{align*}
(here, we have renamed the summation index $p$ as $k$). This solves Exercise
\ref{exe.perm.rearrangement} \textbf{(f)}.
\end{verlong}

\begin{vershort}
\textbf{(g)} Fix $k\in\left[  n\right]  $. Then, $i_{k}\in\left[  k\right]  $
(by (\ref{sol.perm.rearrangement.1})), so that $i_{k}\leq k$ and therefore
$x_{i_{k}}\geq x_{k}$ (since $x_{1}\geq x_{2}\geq\cdots\geq x_{n}$). Hence,
$x_{i_{k}}-x_{k}\geq0$.

Also, $m_{k}\in\left[  k\right]  $ (by Exercise \ref{exe.perm.rearrangement}
\textbf{(b)}), so that $m_{k}\leq k$ and thus $y_{m_{k}}\geq y_{k}$ (since
$y_{1}\geq y_{2}\geq\cdots\geq y_{n}$). Hence, $y_{m_{k}}-y_{k}\geq0$. Now,%
\begin{equation}
\underbrace{\left(  x_{i_{k}}-x_{k}\right)  }_{\geq0}\underbrace{\left(
y_{m_{k}}-y_{k}\right)  }_{\geq0}\geq0. \label{sol.perm.rearrangement.g.1}%
\end{equation}

Now, forget that we fixed $k$. We thus have proven
(\ref{sol.perm.rearrangement.g.1}) for each $k\in\left[  n\right]  $. Now,
Exercise \ref{exe.perm.rearrangement} \textbf{(f)} yields
\[
\sum_{k=1}^{n}x_{k}y_{k}-\sum_{k=1}^{n}x_{k}y_{\sigma\left(  k\right)  }%
=\sum_{k=1}^{n}\underbrace{\left(  x_{i_{k}}-x_{k}\right)  \left(  y_{m_{k}%
}-y_{k}\right)  }_{\substack{\geq0\\\text{(by
(\ref{sol.perm.rearrangement.g.1}))}}}\geq0.
\]
In other words,%
\[
\sum_{k=1}^{n}x_{k}y_{k}\geq\sum_{k=1}^{n}x_{k}y_{\sigma\left(  k\right)  }.
\]
This solves Exercise \ref{exe.perm.rearrangement} \textbf{(g)}. \qedhere

\end{vershort}

\begin{verlong}
\textbf{(g)} Fix $k\in\left\{  1,2,\ldots,n\right\}  $. Then, $k\leq n$ and
thus $\left\{  1,2,\ldots,k\right\}  \subseteq\left\{  1,2,\ldots,n\right\}
$. Hence, $\left[  k\right]  =\left\{  1,2,\ldots,k\right\}  \subseteq\left\{
1,2,\ldots,n\right\}  =\left[  n\right]  $. Also, $k\in\left\{  1,2,\ldots
,n\right\}  =\left[  n\right]  $.

Hence, (\ref{sol.perm.rearrangement.1}) (applied to $j=k$) yields $i_{k}%
\in\left[  k\right]  =\left\{  1,2,\ldots,k\right\}  $, so that $i_{k}\leq k$.
But we have $x_{1}\geq x_{2}\geq\cdots\geq x_{n}$. In other words, if $u$ and
$v$ are two elements of $\left[  n\right]  $ such that $u\leq v$, then
$x_{u}\geq x_{v}$. Applying this to $u=i_{k}$ and $v=k$, we obtain $x_{i_{k}%
}\geq x_{k}$ (since $i_{k}\in\left[  k\right]  \subseteq\left[  n\right]  $
and $k\in\left[  n\right]  $ and $i_{k}\leq k$). Hence, $x_{i_{k}}-x_{k}\geq0$.

Also, Exercise \ref{exe.perm.rearrangement} \textbf{(b)} yields $m_{k}%
\in\left[  k\right]  =\left\{  1,2,\ldots,k\right\}  $ and therefore
$m_{k}\leq k$. But we have $y_{1}\geq y_{2}\geq\cdots\geq y_{n}$. In other
words, if $u$ and $v$ are two elements of $\left[  n\right]  $ such that
$u\leq v$, then $y_{u}\geq y_{v}$. Applying this to $u=m_{k}$ and $v=k$, we
obtain $y_{m_{k}}\geq y_{k}$ (since $m_{k}\in\left[  k\right]  \subseteq
\left[  n\right]  $ and $k\in\left[  n\right]  $ and $m_{k}\leq k$). Hence,
$y_{m_{k}}-y_{k}\geq0$.

Now,%
\begin{equation}
\underbrace{\left(  x_{i_{k}}-x_{k}\right)  }_{\geq0}\underbrace{\left(
y_{m_{k}}-y_{k}\right)  }_{\geq0}\geq0.
\label{sol.perm.rearrangement.g.long.1}%
\end{equation}

Now, forget that we fixed $k$. We thus have proven
(\ref{sol.perm.rearrangement.g.long.1}) for each $k\in\left\{  1,2,\ldots
,n\right\}  $. Now, Exercise \ref{exe.perm.rearrangement} \textbf{(f)} yields
\[
\sum_{k=1}^{n}x_{k}y_{k}-\sum_{k=1}^{n}x_{k}y_{\sigma\left(  k\right)  }%
=\sum_{k=1}^{n}\underbrace{\left(  x_{i_{k}}-x_{k}\right)  \left(  y_{m_{k}%
}-y_{k}\right)  }_{\substack{\geq0\\\text{(by
(\ref{sol.perm.rearrangement.g.long.1}))}}}\geq0.
\]
In other words,%
\[
\sum_{k=1}^{n}x_{k}y_{k}\geq\sum_{k=1}^{n}x_{k}y_{\sigma\left(  k\right)  }.
\]
This solves Exercise \ref{exe.perm.rearrangement} \textbf{(g)}.
\end{verlong}
\end{proof}

\subsection{\label{sect.sol.perm.sigmacrosstau}Solution to Exercise
\ref{exe.perm.sigmacrosstau}}

Before we come to the solution of Exercise \ref{exe.perm.sigmacrosstau}, we
need to lay a lot of groundwork. First, we introduce a notation, which
generalizes the definition of $\sigma\times\tau$ in Exercise
\ref{exe.perm.sigmacrosstau}:

\begin{definition}
\label{def.sol.sigmacrosstau.aXb}If $X$, $X^{\prime}$, $Y$ and $Y^{\prime}$
are four sets and if $\alpha:X\rightarrow X^{\prime}$ and $\beta:Y\rightarrow
Y^{\prime}$ are two maps, then $\alpha\times\beta$ will denote the map%
\begin{align*}
X\times Y  &  \rightarrow X^{\prime}\times Y^{\prime},\\
\left(  x,y\right)   &  \mapsto\left(  \alpha\left(  x\right)  ,\beta\left(
y\right)  \right)  .
\end{align*}

\end{definition}

We shall use Definition \ref{def.sol.sigmacrosstau.aXb} throughout Section
\ref{sect.sol.perm.sigmacrosstau}.

The following properties of this definition are straightforward to prove:

\begin{proposition}
\label{prop.sol.sigmacrosstau.fgf'g'}Let $X$, $X^{\prime}$, $X^{\prime\prime}%
$, $Y$, $Y^{\prime}$ and $Y^{\prime\prime}$ be six sets. Let $\alpha
:X\rightarrow X^{\prime}$, $\alpha^{\prime}:X^{\prime}\rightarrow
X^{\prime\prime}$, $\beta:Y\rightarrow Y^{\prime}$ and $\beta^{\prime
}:Y^{\prime}\rightarrow Y^{\prime\prime}$ be four maps. Then,%
\[
\left(  \alpha^{\prime}\times\beta^{\prime}\right)  \circ\left(  \alpha
\times\beta\right)  =\left(  \alpha^{\prime}\circ\alpha\right)  \times\left(
\beta^{\prime}\circ\beta\right)  .
\]

\end{proposition}

\begin{vershort}
\begin{proof}
[Proof of Proposition \ref{prop.sol.sigmacrosstau.fgf'g'}.]Straightforward
computation reveals that the maps $\left(  \alpha^{\prime}\times\beta^{\prime
}\right)  \circ\left(  \alpha\times\beta\right)  $ and $\left(  \alpha
^{\prime}\circ\alpha\right)  \times\left(  \beta^{\prime}\circ\beta\right)  $
send any given element $\left(  x,y\right)  \in X\times Y$ to the same image
(namely, to $\left(  \alpha^{\prime}\left(  \alpha\left(  x\right)  \right)
,\beta^{\prime}\left(  \beta\left(  y\right)  \right)  \right)  $). Thus,
these two maps are identical. This proves Proposition
\ref{prop.sol.sigmacrosstau.fgf'g'}.
\end{proof}
\end{vershort}

\begin{verlong}
\begin{proof}
[Proof of Proposition \ref{prop.sol.sigmacrosstau.fgf'g'}.]Let $u\in X\times
Y$ be arbitrary. Thus, $u=\left(  x,y\right)  $ for some $x\in X$ and $y\in
Y$. Consider these $x$ and $y$. Now,%
\[
\left(  \alpha\times\beta\right)  \left(  \underbrace{u}_{=\left(  x,y\right)
}\right)  =\left(  \alpha\times\beta\right)  \left(  \left(  x,y\right)
\right)  =\left(  \alpha\left(  x\right)  ,\beta\left(  y\right)  \right)
\]
(by the definition of $\alpha\times\beta$). Comparing the equality%
\begin{align*}
\left(  \left(  \alpha^{\prime}\times\beta^{\prime}\right)  \circ\left(
\alpha\times\beta\right)  \right)  \left(  u\right)   &  =\left(
\alpha^{\prime}\times\beta^{\prime}\right)  \left(  \underbrace{\left(
\alpha\times\beta\right)  \left(  u\right)  }_{=\left(  \alpha\left(
x\right)  ,\beta\left(  y\right)  \right)  }\right)  =\left(  \alpha^{\prime
}\times\beta^{\prime}\right)  \left(  \alpha\left(  x\right)  ,\beta\left(
y\right)  \right) \\
&  =\left(  \underbrace{\alpha^{\prime}\left(  \alpha\left(  x\right)
\right)  }_{=\left(  \alpha^{\prime}\circ\alpha\right)  \left(  x\right)
},\underbrace{\beta^{\prime}\left(  \beta\left(  y\right)  \right)
}_{=\left(  \beta^{\prime}\circ\beta\right)  \left(  y\right)  }\right) \\
&  \ \ \ \ \ \ \ \ \ \ \left(  \text{by the definition of }\alpha^{\prime
}\times\beta^{\prime}\right) \\
&  =\left(  \left(  \alpha^{\prime}\circ\alpha\right)  \left(  x\right)
,\left(  \beta^{\prime}\circ\beta\right)  \left(  y\right)  \right)
\end{align*}
with the equality%
\begin{align*}
\left(  \left(  \alpha^{\prime}\circ\alpha\right)  \times\left(  \beta
^{\prime}\circ\beta\right)  \right)  \left(  \underbrace{u}_{=\left(
x,y\right)  }\right)   &  =\left(  \left(  \alpha^{\prime}\circ\alpha\right)
\times\left(  \beta^{\prime}\circ\beta\right)  \right)  \left(  \left(
x,y\right)  \right)  =\left(  \left(  \alpha^{\prime}\circ\alpha\right)
\left(  x\right)  ,\left(  \beta^{\prime}\circ\beta\right)  \left(  y\right)
\right) \\
&  \ \ \ \ \ \ \ \ \ \ \left(  \text{by the definition of }\left(
\alpha^{\prime}\circ\alpha\right)  \times\left(  \beta^{\prime}\circ
\beta\right)  \right)  ,
\end{align*}
we obtain $\left(  \left(  \alpha^{\prime}\times\beta^{\prime}\right)
\circ\left(  \alpha\times\beta\right)  \right)  \left(  u\right)  =\left(
\left(  \alpha^{\prime}\circ\alpha\right)  \times\left(  \beta^{\prime}%
\circ\beta\right)  \right)  \left(  u\right)  $.

Now, forget that we fixed $u$. We thus have shown that $\left(  \left(
\alpha^{\prime}\times\beta^{\prime}\right)  \circ\left(  \alpha\times
\beta\right)  \right)  \left(  u\right)  =\left(  \left(  \alpha^{\prime}%
\circ\alpha\right)  \times\left(  \beta^{\prime}\circ\beta\right)  \right)
\left(  u\right)  $ for each $u\in X\times Y$. In other words, $\left(
\alpha^{\prime}\times\beta^{\prime}\right)  \circ\left(  \alpha\times
\beta\right)  =\left(  \alpha^{\prime}\circ\alpha\right)  \times\left(
\beta^{\prime}\circ\beta\right)  $. This proves Proposition
\ref{prop.sol.sigmacrosstau.fgf'g'}.
\end{proof}
\end{verlong}

\begin{proposition}
\label{prop.sol.sigmacrosstau.idxid}Let $U$ and $V$ be two sets. Then,
$\operatorname*{id}\nolimits_{U}\times\operatorname*{id}\nolimits_{V}%
=\operatorname*{id}\nolimits_{U\times V}$.
\end{proposition}

\begin{vershort}
\begin{proof}
[Proof of Proposition \ref{prop.sol.sigmacrosstau.idxid}.]Again, this follows
from the straightforward computation of $\left(  \operatorname*{id}%
\nolimits_{U}\times\operatorname*{id}\nolimits_{V}\right)  \left(  x,y\right)
$ for each $\left(  x,y\right)  \in U\times V$.
\end{proof}
\end{vershort}

\begin{verlong}
\begin{proof}
[Proof of Proposition \ref{prop.sol.sigmacrosstau.idxid}.]Let $c\in U\times V$
be arbitrary. Thus, $c=\left(  u,v\right)  $ for some $u\in U$ and $v\in V$.
Consider these $u$ and $v$. Now,%
\begin{align*}
\left(  \operatorname*{id}\nolimits_{U}\times\operatorname*{id}\nolimits_{V}%
\right)  \left(  \underbrace{c}_{=\left(  u,v\right)  }\right)   &  =\left(
\operatorname*{id}\nolimits_{U}\times\operatorname*{id}\nolimits_{V}\right)
\left(  \left(  u,v\right)  \right)  =\left(  \underbrace{\operatorname*{id}%
\nolimits_{U}\left(  u\right)  }_{=u},\underbrace{\operatorname*{id}%
\nolimits_{V}\left(  v\right)  }_{=v}\right) \\
&  =\left(  u,v\right)  =c=\operatorname*{id}\nolimits_{U\times V}\left(
c\right)  .
\end{align*}

Now, forget that we fixed $c$. We thus have proven that $\left(
\operatorname*{id}\nolimits_{U}\times\operatorname*{id}\nolimits_{V}\right)
\left(  c\right)  =\operatorname*{id}\nolimits_{U\times V}\left(  c\right)  $
for each $c\in U\times V$. In other words, $\operatorname*{id}\nolimits_{U}%
\times\operatorname*{id}\nolimits_{V}=\operatorname*{id}\nolimits_{U\times V}%
$. This proves Proposition \ref{prop.sol.sigmacrosstau.idxid}.
\end{proof}
\end{verlong}

\begin{corollary}
\label{cor.sol.sigmacrosstau.fixid}Let $U$ and $V$ be two sets. Let
$k\in\mathbb{N}$. Let $f_{1},f_{2},\ldots,f_{k}$ be $k$ maps from $U$ to $U$.
For each $i\in\left\{  1,2,\ldots,k\right\}  $, define a map $g_{i}:U\times
V\rightarrow U\times V$ by $g_{i}=f_{i}\times\operatorname*{id}\nolimits_{V}$.
Then,
\[
g_{1}\circ g_{2}\circ\cdots\circ g_{k}=\left(  f_{1}\circ f_{2}\circ
\cdots\circ f_{k}\right)  \times\operatorname*{id}\nolimits_{V}.
\]

\end{corollary}

\begin{proof}
[Proof of Corollary \ref{cor.sol.sigmacrosstau.fixid}.]We claim that%
\begin{equation}
g_{1}\circ g_{2}\circ\cdots\circ g_{m}=\left(  f_{1}\circ f_{2}\circ
\cdots\circ f_{m}\right)  \times\operatorname*{id}\nolimits_{V}
\label{pf.cor.sol.sigmacrosstau.fixid.goal}%
\end{equation}
for each $m\in\left\{  0,1,\ldots,k\right\}  $.

[\textit{Proof of (\ref{pf.cor.sol.sigmacrosstau.fixid.goal}):} We shall prove
(\ref{pf.cor.sol.sigmacrosstau.fixid.goal}) by induction on $m$:

\textit{Induction base:} We have%
\[
g_{1}\circ g_{2}\circ\cdots\circ g_{0}=\left(  \text{empty composition of maps
}U\times V\rightarrow U\times V\right)  =\operatorname*{id}\nolimits_{U\times
V}.
\]
Comparing this with%
\[
\underbrace{\left(  f_{1}\circ f_{2}\circ\cdots\circ f_{0}\right)
}_{\substack{=\left(  \text{empty composition of maps }U\rightarrow U\right)
\\=\operatorname*{id}\nolimits_{U}}}\times\operatorname*{id}\nolimits_{V}%
=\operatorname*{id}\nolimits_{U}\times\operatorname*{id}\nolimits_{V}%
=\operatorname*{id}\nolimits_{U\times V}\ \ \ \ \ \ \ \ \ \ \left(  \text{by
Proposition \ref{prop.sol.sigmacrosstau.idxid}}\right)  ,
\]
we obtain $g_{1}\circ g_{2}\circ\cdots\circ g_{0}=\left(  f_{1}\circ
f_{2}\circ\cdots\circ f_{0}\right)  \times\operatorname*{id}\nolimits_{V}$. In
other words, (\ref{pf.cor.sol.sigmacrosstau.fixid.goal}) holds for $m=0$. This
completes the induction base.

\textit{Induction step:} Let $M\in\left\{  0,1,\ldots,k\right\}  $ be
positive. Assume that (\ref{pf.cor.sol.sigmacrosstau.fixid.goal}) holds for
$m=M-1$. We must prove that (\ref{pf.cor.sol.sigmacrosstau.fixid.goal}) holds
for $m=M$.

\begin{verlong}
We have $M\neq0$ (since $M$ is positive). Combining $M\in\left\{
0,1,\ldots,k\right\}  $ with $M\neq0$, we obtain $M\in\left\{  0,1,\ldots
,k\right\}  \setminus\left\{  0\right\}  =\left\{  1,2,\ldots,k\right\}  $.
Thus, $M-1\in\left\{  0,1,\ldots,k-1\right\}  \subseteq\left\{  0,1,\ldots
,k\right\}  $.
\end{verlong}

We have assumed that (\ref{pf.cor.sol.sigmacrosstau.fixid.goal}) holds for
$m=M-1$. In other words, we have%
\[
g_{1}\circ g_{2}\circ\cdots\circ g_{M-1}=\left(  f_{1}\circ f_{2}\circ
\cdots\circ f_{M-1}\right)  \times\operatorname*{id}\nolimits_{V}.
\]

But $M$ is positive; thus,%
\begin{align*}
&  g_{1}\circ g_{2}\circ\cdots\circ g_{M}\\
&  =\underbrace{\left(  g_{1}\circ g_{2}\circ\cdots\circ g_{M-1}\right)
}_{=\left(  f_{1}\circ f_{2}\circ\cdots\circ f_{M-1}\right)  \times
\operatorname*{id}\nolimits_{V}}\circ\underbrace{g_{M}}_{\substack{=f_{M}%
\times\operatorname*{id}\nolimits_{V}\\\text{(by the definition of }%
g_{M}\text{)}}}\\
&  =\left(  \left(  f_{1}\circ f_{2}\circ\cdots\circ f_{M-1}\right)
\times\operatorname*{id}\nolimits_{V}\right)  \circ\left(  f_{M}%
\times\operatorname*{id}\nolimits_{V}\right) \\
&  =\underbrace{\left(  \left(  f_{1}\circ f_{2}\circ\cdots\circ
f_{M-1}\right)  \circ f_{M}\right)  }_{=f_{1}\circ f_{2}\circ\cdots\circ
f_{M}}\times\underbrace{\left(  \operatorname*{id}\nolimits_{V}\circ
\operatorname*{id}\nolimits_{V}\right)  }_{=\operatorname*{id}\nolimits_{V}}\\
&  \ \ \ \ \ \ \ \ \ \ \left(
\begin{array}
[c]{c}%
\text{by Proposition \ref{prop.sol.sigmacrosstau.fgf'g'}}\\
\text{(applied to }X=U\text{, }X^{\prime}=U\text{, }X^{\prime\prime}=U\text{,
}Y=V\text{, }Y^{\prime}=V\text{, }Y^{\prime\prime}=V\text{,}\\
\alpha=f_{M}\text{, }\alpha^{\prime}=f_{1}\circ f_{2}\circ\cdots\circ
f_{M-1}\text{, }\beta=\operatorname*{id}\nolimits_{V}\text{ and }\beta
^{\prime}=\operatorname*{id}\nolimits_{V}\text{)}%
\end{array}
\right) \\
&  =\left(  f_{1}\circ f_{2}\circ\cdots\circ f_{M}\right)  \times
\operatorname*{id}\nolimits_{V}.
\end{align*}
In other words, (\ref{pf.cor.sol.sigmacrosstau.fixid.goal}) holds for $m=M$.
This completes the induction step. Thus,
(\ref{pf.cor.sol.sigmacrosstau.fixid.goal}) is proven.]

Now, $k\in\left\{  0,1,\ldots,k\right\}  $ (since $k\in\mathbb{N}$). Hence,
(\ref{pf.cor.sol.sigmacrosstau.fixid.goal}) (applied to $m=k$) shows that
$g_{1}\circ g_{2}\circ\cdots\circ g_{k}=\left(  f_{1}\circ f_{2}\circ
\cdots\circ f_{k}\right)  \times\operatorname*{id}\nolimits_{V}$. This proves
Corollary \ref{cor.sol.sigmacrosstau.fixid}.
\end{proof}

\begin{corollary}
\label{cor.sol.sigmacrosstau.bijxbij}Let $X$, $X^{\prime}$, $Y$ and
$Y^{\prime}$ be four sets. Let $\alpha:X\rightarrow X^{\prime}$ and
$\beta:Y\rightarrow Y^{\prime}$ be two bijective maps. Then, the map
$\alpha\times\beta:X\times Y\rightarrow X^{\prime}\times Y^{\prime}$ is
bijective as well, and its inverse is the map $\alpha^{-1}\times\beta
^{-1}:X^{\prime}\times Y^{\prime}\rightarrow X\times Y$.
\end{corollary}

\begin{vershort}
\begin{proof}
[Proof of Corollary \ref{cor.sol.sigmacrosstau.bijxbij}.]This can be derived
from Proposition \ref{prop.sol.sigmacrosstau.idxid} and Proposition
\ref{prop.sol.sigmacrosstau.fgf'g'} (or, again, checked by straightforward computation).
\end{proof}
\end{vershort}

\begin{verlong}
\begin{proof}
[Proof of Corollary \ref{cor.sol.sigmacrosstau.bijxbij}.]Proposition
\ref{prop.sol.sigmacrosstau.idxid} (applied to $U=X$ and $V=Y$) yields
$\operatorname*{id}\nolimits_{X}\times\operatorname*{id}\nolimits_{Y}%
=\operatorname*{id}\nolimits_{X\times Y}$. Proposition
\ref{prop.sol.sigmacrosstau.idxid} (applied to $U=X^{\prime}$ and
$V=Y^{\prime}$) yields $\operatorname*{id}\nolimits_{X^{\prime}}%
\times\operatorname*{id}\nolimits_{Y^{\prime}}=\operatorname*{id}%
\nolimits_{X^{\prime}\times Y^{\prime}}$.

The map $\alpha$ is bijective. Thus, its inverse map $\alpha^{-1}:X^{\prime
}\rightarrow X$ is well-defined.

The map $\beta$ is bijective. Thus, its inverse map $\beta^{-1}:Y^{\prime
}\rightarrow Y$ is well-defined.

Thus, the map $\alpha^{-1}\times\beta^{-1}:X^{\prime}\times Y^{\prime
}\rightarrow X\times Y$ is well-defined. Furthermore, Proposition
\ref{prop.sol.sigmacrosstau.fgf'g'} (applied to $X^{\prime\prime}=X$,
$Y^{\prime\prime}=Y$, $\alpha^{\prime}=\alpha^{-1}$ and $\beta^{\prime}%
=\beta^{-1}$) yields%
\[
\left(  \alpha^{-1}\times\beta^{-1}\right)  \circ\left(  \alpha\times
\beta\right)  =\underbrace{\left(  \alpha^{-1}\circ\alpha\right)
}_{=\operatorname*{id}\nolimits_{X}}\times\underbrace{\left(  \beta^{-1}%
\circ\beta\right)  }_{=\operatorname*{id}\nolimits_{Y}}=\operatorname*{id}%
\nolimits_{X}\times\operatorname*{id}\nolimits_{Y}=\operatorname*{id}%
\nolimits_{X\times Y}.
\]
Also, Proposition \ref{prop.sol.sigmacrosstau.fgf'g'} (applied to $X^{\prime}%
$, $X$, $X^{\prime}$, $Y^{\prime}$, $Y$, $Y^{\prime}$, $\alpha^{-1}$, $\alpha
$, $\beta^{-1}$ and $\beta$ instead of $X$, $X^{\prime}$, $X^{\prime\prime}$,
$Y$, $Y^{\prime}$, $Y^{\prime\prime}$, $\alpha$, $\alpha^{\prime}$, $\beta$
and $\beta^{\prime}$) yields%
\[
\left(  \alpha\times\beta\right)  \circ\left(  \alpha^{-1}\times\beta
^{-1}\right)  =\underbrace{\left(  \alpha\circ\alpha^{-1}\right)
}_{=\operatorname*{id}\nolimits_{X^{\prime}}}\times\underbrace{\left(
\beta\circ\beta^{-1}\right)  }_{=\operatorname*{id}\nolimits_{Y^{\prime}}%
}=\operatorname*{id}\nolimits_{X^{\prime}}\times\operatorname*{id}%
\nolimits_{Y^{\prime}}=\operatorname*{id}\nolimits_{X^{\prime}\times
Y^{\prime}}.
\]

The maps $\alpha\times\beta$ and $\alpha^{-1}\times\beta^{-1}$ are mutually
inverse (since $\left(  \alpha^{-1}\times\beta^{-1}\right)  \circ\left(
\alpha\times\beta\right)  =\operatorname*{id}\nolimits_{X\times Y}$ and
$\left(  \alpha\times\beta\right)  \circ\left(  \alpha^{-1}\times\beta
^{-1}\right)  =\operatorname*{id}\nolimits_{X^{\prime}\times Y^{\prime}}$).
Thus, the map $\alpha\times\beta$ is invertible, and hence is bijective. So we
have proven that the map $\alpha\times\beta:X\times Y\rightarrow X^{\prime
}\times Y^{\prime}$ is bijective. Furthermore, its inverse is the map
$\alpha^{-1}\times\beta^{-1}:X^{\prime}\times Y^{\prime}\rightarrow X\times Y$
(since the maps $\alpha\times\beta$ and $\alpha^{-1}\times\beta^{-1}$ are
mutually inverse). This completes the proof of Corollary
\ref{cor.sol.sigmacrosstau.bijxbij}.
\end{proof}
\end{verlong}

\begin{corollary}
\label{cor.sol.sigmacrosstau.permxperm}Let $U$ and $V$ be two sets. Let
$\sigma$ be a permutation of $U$. Let $\tau$ be a permutation of $V$. Then,
$\sigma\times\tau$ is a permutation of $U\times V$.
\end{corollary}

\begin{vershort}
\begin{proof}
[Proof of Corollary \ref{cor.sol.sigmacrosstau.permxperm}.]The map $\sigma$ is
a permutation of $U$. In other words, $\sigma$ is a bijective map
$U\rightarrow U$. Similarly, $\tau$ is a bijective map $V\rightarrow V$.
Hence, Corollary \ref{cor.sol.sigmacrosstau.bijxbij} (applied to $X=U$, $Y=U$,
$Y=V$, $Y^{\prime}=V$, $\alpha=\sigma$ and $\beta=\tau$) shows that the map
$\sigma\times\tau:U\times V\rightarrow U\times V$ is bijective as well. In
other words, $\sigma\times\tau$ is a permutation of $U\times V$. This proves
Corollary \ref{cor.sol.sigmacrosstau.permxperm}.
\end{proof}
\end{vershort}

\begin{verlong}
\begin{proof}
[Proof of Corollary \ref{cor.sol.sigmacrosstau.permxperm}.]We know that
$\sigma$ is a permutation of $U$. In other words, $\sigma$ is a bijection from
$U$ to $U$. In other words, $\sigma$ is a bijective map $U\rightarrow U$.

We know that $\tau$ is a permutation of $V$. In other words, $\tau$ is a
bijection from $V$ to $V$. In other words, $\tau$ is a bijective map
$V\rightarrow V$.

Corollary \ref{cor.sol.sigmacrosstau.bijxbij} (applied to $X=U$, $Y=U$, $Y=V$,
$Y^{\prime}=V$, $\alpha=\sigma$ and $\beta=\tau$) shows that the map
$\sigma\times\tau:U\times V\rightarrow U\times V$ is bijective as well, and
its inverse is the map $\sigma^{-1}\times\tau^{-1}:U\times V\rightarrow
U\times V$. In particular, $\sigma\times\tau$ is a bijective map $U\times
V\rightarrow U\times V$. In other words, $\sigma\times\tau$ is a bijection
from $U\times V$ to $U\times V$. In other words, $\sigma\times\tau$ is a
permutation of $U\times V$. This proves Corollary
\ref{cor.sol.sigmacrosstau.permxperm}.
\end{proof}
\end{verlong}

\begin{proposition}
\label{prop.sol.sigmacrosstau.trans}Let $U$ and $V$ be two sets. Let
$\gamma:U\times V\rightarrow V\times U$ be the map defined by%
\[
\left(  \gamma\left(  \left(  u,v\right)  \right)  =\left(  v,u\right)
\ \ \ \ \ \ \ \ \ \ \text{for each }\left(  u,v\right)  \in U\times V\right)
.
\]
Then:

\textbf{(a)} The map $\gamma:U\times V\rightarrow V\times U$ is a bijection.

\textbf{(b)} Let $f$ be any map $U\rightarrow U$. Let $g$ be any map
$V\rightarrow V$. Then, $f\times g=\gamma^{-1}\circ\left(  g\times f\right)
\circ\gamma$.
\end{proposition}

\begin{vershort}
\begin{proof}
[Proof of Proposition \ref{prop.sol.sigmacrosstau.trans}.]\textbf{(a)} Let
$\delta:V\times U\rightarrow U\times V$ be the map defined by%
\[
\left(  \delta\left(  \left(  v,u\right)  \right)  =\left(  u,v\right)
\ \ \ \ \ \ \ \ \ \ \text{for each }\left(  v,u\right)  \in V\times U\right)
.
\]
Clearly, the maps $\gamma$ and $\delta$ are mutually inverse. Thus, the map
$\gamma$ is invertible, i.e., is a bijection. This proves Proposition
\ref{prop.sol.sigmacrosstau.trans} \textbf{(a)}.

\textbf{(b)} For each $\left(  u,v\right)  \in U\times V$, we have
\begin{align*}
&  \left(  \gamma\circ\left(  f\times g\right)  \right)  \left(  \left(
u,v\right)  \right) \\
&  =\gamma\left(  \underbrace{\left(  f\times g\right)  \left(  \left(
u,v\right)  \right)  }_{\substack{=\left(  f\left(  u\right)  ,g\left(
v\right)  \right)  \\\text{(by the definition of }f\times g\text{)}}}\right)
\\
&  =\gamma\left(  \left(  f\left(  u\right)  ,g\left(  v\right)  \right)
\right)  =\left(  g\left(  v\right)  ,f\left(  u\right)  \right)
\ \ \ \ \ \ \ \ \ \ \left(  \text{by the definition of }\gamma\right)
\end{align*}
and%
\begin{align*}
&  \left(  \left(  g\times f\right)  \circ\gamma\right)  \left(  \left(
u,v\right)  \right) \\
&  =\left(  g\times f\right)  \left(  \underbrace{\gamma\left(  \left(
u,v\right)  \right)  }_{=\left(  v,u\right)  }\right) \\
&  =\left(  g\times f\right)  \left(  \left(  v,u\right)  \right)  =\left(
g\left(  v\right)  ,f\left(  u\right)  \right)  \ \ \ \ \ \ \ \ \ \ \left(
\text{by the definition of }g\times f\right)  .
\end{align*}
Comparing these two equalities, we conclude that
\[
\left(  \gamma\circ\left(  f\times g\right)  \right)  \left(  \left(
u,v\right)  \right)  =\left(  \left(  g\times f\right)  \circ\gamma\right)
\left(  \left(  u,v\right)  \right)  \ \ \ \ \ \ \ \ \ \ \text{for each
}\left(  u,v\right)  \in U\times V.
\]
In other words, $\gamma\circ\left(  f\times g\right)  =\left(  g\times
f\right)  \circ\gamma$. Thus,%
\[
\gamma^{-1}\circ\underbrace{\left(  g\times f\right)  \circ\gamma}%
_{=\gamma\circ\left(  f\times g\right)  }=\underbrace{\gamma^{-1}\circ\gamma
}_{=\operatorname*{id}\nolimits_{U\times V}}\circ\left(  f\times g\right)
=f\times g.
\]
This proves Proposition \ref{prop.sol.sigmacrosstau.trans} \textbf{(b)}.
\end{proof}
\end{vershort}

\begin{verlong}
\begin{proof}
[Proof of Proposition \ref{prop.sol.sigmacrosstau.trans}.]\textbf{(a)} Let
$\delta:V\times U\rightarrow U\times V$ be the map defined by%
\[
\left(  \delta\left(  \left(  v,u\right)  \right)  =\left(  u,v\right)
\ \ \ \ \ \ \ \ \ \ \text{for each }\left(  v,u\right)  \in V\times U\right)
.
\]
Then, $\delta\circ\gamma=\operatorname*{id}\nolimits_{U\times V}%
$\ \ \ \ \footnote{\textit{Proof.} Let $c\in U\times V$. Thus, $c$ can be
written in the form $c=\left(  u,v\right)  $ for some $u\in U$ and $v\in V$.
Consider these $u$ and $v$.
\par
We have
\begin{align*}
\left(  \delta\circ\gamma\right)  \left(  \underbrace{c}_{=\left(  u,v\right)
}\right)   &  =\left(  \delta\circ\gamma\right)  \left(  \left(  u,v\right)
\right)  =\delta\left(  \underbrace{\gamma\left(  \left(  u,v\right)  \right)
}_{\substack{=\left(  v,u\right)  \\\text{(by the definition of }%
\gamma\text{)}}}\right)  =\delta\left(  \left(  v,u\right)  \right)  =\left(
u,v\right) \\
&  =c=\operatorname*{id}\nolimits_{U\times V}\left(  c\right)  .
\end{align*}
\par
Now, forget that we fixed $c$. Thus, we have shown that $\left(  \delta
\circ\gamma\right)  \left(  c\right)  =\operatorname*{id}\nolimits_{U\times
V}\left(  c\right)  $ for each $c\in U\times V$. In other words, $\delta
\circ\gamma=\operatorname*{id}\nolimits_{U\times V}$.}. The same argument
(applied to $V$, $U$, $\delta$ and $\gamma$ instead of $U$, $V$, $\gamma$ and
$\delta$) shows that $\gamma\circ\delta=\operatorname*{id}\nolimits_{V\times
U}$. Now, the maps $\gamma$ and $\delta$ are mutually inverse (since
$\delta\circ\gamma=\operatorname*{id}\nolimits_{U\times V}$ and $\gamma
\circ\delta=\operatorname*{id}\nolimits_{V\times U}$). Hence, the map $\gamma$
is invertible, and therefore bijective. In other words, the map $\gamma
:U\times V\rightarrow V\times U$ is a bijection. This proves Proposition
\ref{prop.sol.sigmacrosstau.trans} \textbf{(a)}.

\textbf{(b)} Let $c\in U\times V$. Thus, $c$ can be written in the form
$c=\left(  u,v\right)  $ for some $u\in U$ and $v\in V$. Consider these $u$
and $v$.

We have%
\begin{align*}
\left(  \gamma\circ\left(  f\times g\right)  \right)  \left(  \underbrace{c}%
_{=\left(  u,v\right)  }\right)   &  =\left(  \gamma\circ\left(  f\times
g\right)  \right)  \left(  \left(  u,v\right)  \right)  =\gamma\left(
\underbrace{\left(  f\times g\right)  \left(  \left(  u,v\right)  \right)
}_{\substack{=\left(  f\left(  u\right)  ,g\left(  v\right)  \right)
\\\text{(by the definition of }f\times g\text{)}}}\right) \\
&  =\gamma\left(  \left(  f\left(  u\right)  ,g\left(  v\right)  \right)
\right)  =\left(  g\left(  v\right)  ,f\left(  u\right)  \right)
\end{align*}
(by the definition of $\gamma$). Comparing this with%
\begin{align*}
\left(  \left(  g\times f\right)  \circ\gamma\right)  \left(  \underbrace{c}%
_{=\left(  u,v\right)  }\right)   &  =\left(  \left(  g\times f\right)
\circ\gamma\right)  \left(  \left(  u,v\right)  \right)  =\left(  g\times
f\right)  \left(  \underbrace{\gamma\left(  \left(  u,v\right)  \right)
}_{\substack{=\left(  v,u\right)  \\\text{(by the definition of }%
\gamma\text{)}}}\right) \\
&  =\left(  g\times f\right)  \left(  \left(  v,u\right)  \right)  =\left(
g\left(  v\right)  ,f\left(  u\right)  \right)  \ \ \ \ \ \ \ \ \ \ \left(
\text{by the definition of }g\times f\right)  ,
\end{align*}
we obtain $\left(  \gamma\circ\left(  f\times g\right)  \right)  \left(
c\right)  =\left(  \left(  g\times f\right)  \circ\gamma\right)  \left(
c\right)  $.

Now, forget that we fixed $c$. We thus have proven that $\left(  \gamma
\circ\left(  f\times g\right)  \right)  \left(  c\right)  =\left(  \left(
g\times f\right)  \circ\gamma\right)  \left(  c\right)  $ for each $c\in
U\times V$. In other words, $\gamma\circ\left(  f\times g\right)  =\left(
g\times f\right)  \circ\gamma$. Thus, $\left(  g\times f\right)  \circ
\gamma=\gamma\circ\left(  f\times g\right)  $, so that%
\[
\gamma^{-1}\circ\underbrace{\left(  g\times f\right)  \circ\gamma}%
_{=\gamma\circ\left(  f\times g\right)  }=\underbrace{\gamma^{-1}\circ\gamma
}_{=\operatorname*{id}\nolimits_{U\times V}}\circ\left(  f\times g\right)
=f\times g.
\]
In other words, $f\times g=\gamma^{-1}\circ\left(  g\times f\right)
\circ\gamma$. This proves Proposition \ref{prop.sol.sigmacrosstau.trans}
\textbf{(b)}.
\end{proof}
\end{verlong}

Next, let us discuss some properties of signs of permutations of finite sets:

\begin{proposition}
\label{prop.sol.sigmacrosstau.signs-conj}Let $X$ and $Y$ be two finite sets.
Let $f:X\rightarrow Y$ be a bijection. Let $\sigma$ be a permutation of $X$.
Then, $f\circ\sigma\circ f^{-1}$ is a permutation of $Y$ and satisfies
$\left(  -1\right)  ^{f\circ\sigma\circ f^{-1}}=\left(  -1\right)  ^{\sigma}$.
\end{proposition}

\begin{proof}
[Proof of Proposition \ref{prop.sol.sigmacrosstau.signs-conj}.]The map
$\sigma$ is a permutation of $X$. In other words, $\sigma$ is a bijection from
$X$ to $X$. In other words, $\sigma$ is a bijective map $X\rightarrow X$.

\begin{vershort}
Also, $f:X\rightarrow Y$ is a bijection. Hence, the inverse $f^{-1}$ of $f$ is
well-defined and is a bijection $Y\rightarrow X$.
\end{vershort}

\begin{verlong}
Also, $f:X\rightarrow Y$ is a bijection. In other words, $f$ is a bijective
map $X\rightarrow Y$. Hence, the inverse $f^{-1}$ of $f$ is well-defined and
is a bijective map $Y\rightarrow X$.
\end{verlong}

\begin{vershort}
The map $f\circ\sigma\circ f^{-1}:Y\rightarrow Y$ is a bijection (since it is
the composition of the bijections $f^{-1}$, $\sigma$ and $f$). In other words,
$f\circ\sigma\circ f^{-1}$ is a permutation of $Y$. Hence, $\left(  -1\right)
^{f\circ\sigma\circ f^{-1}}$ is well-defined.
\end{vershort}

\begin{verlong}
The map $f\circ\sigma\circ f^{-1}:Y\rightarrow Y$ is bijective (since it is
the composition of the three bijective maps $f^{-1}$, $\sigma$ and $f$). In
other words, $f\circ\sigma\circ f^{-1}$ is a bijection from $Y$ to $Y$. In
other words, $f\circ\sigma\circ f^{-1}$ is a permutation of $Y$. Hence,
$\left(  -1\right)  ^{f\circ\sigma\circ f^{-1}}$ is well-defined.
\end{verlong}

It remains to prove that $\left(  -1\right)  ^{f\circ\sigma\circ f^{-1}%
}=\left(  -1\right)  ^{\sigma}$.

Define $n\in\mathbb{N}$ by $n=\left\vert X\right\vert $. (This is
well-defined, since the set $X$ is finite.)

We have $\left\vert X\right\vert =n=\left\vert \left\{  1,2,\ldots,n\right\}
\right\vert $. Hence, there exists a bijection $\psi:X\rightarrow\left\{
1,2,\ldots,n\right\}  $. Consider such a $\psi$.

Recall the definition of $\left(  -1\right)  _{\psi}^{\sigma}$ given in
Exercise \ref{exe.ps4.2}. This definition shows that $\left(  -1\right)
_{\psi}^{\sigma}=\left(  -1\right)  ^{\psi\circ\sigma\circ\psi^{-1}}$.

\begin{vershort}
The map $\psi\circ f^{-1}$ is a bijection from $Y$ to $\left\{  1,2,\ldots
,n\right\}  $ (since $\psi$ and $f^{-1}$ are bijections). Hence, the
definition of $\left(  -1\right)  _{\psi\circ f^{-1}}^{f\circ\sigma\circ
f^{-1}}$ (given in Exercise \ref{exe.ps4.2}) shows that
\[
\left(  -1\right)  _{\psi\circ f^{-1}}^{f\circ\sigma\circ f^{-1}}=\left(
-1\right)  ^{\left(  \psi\circ f^{-1}\right)  \circ\left(  f\circ\sigma\circ
f^{-1}\right)  \circ\left(  \psi\circ f^{-1}\right)  ^{-1}}.
\]
In view of
\begin{align*}
&  \underbrace{\left(  \psi\circ f^{-1}\right)  \circ\left(  f\circ\sigma\circ
f^{-1}\right)  }_{=\psi\circ f^{-1}\circ f\circ\sigma\circ f^{-1}}%
\circ\underbrace{\left(  \psi\circ f^{-1}\right)  ^{-1}}_{=\left(
f^{-1}\right)  ^{-1}\circ\psi^{-1}}\\
&  =\psi\circ\underbrace{f^{-1}\circ f}_{=\operatorname*{id}\nolimits_{X}%
}\circ\sigma\circ\underbrace{f^{-1}\circ\left(  f^{-1}\right)  ^{-1}%
}_{=\operatorname*{id}\nolimits_{X}}\circ\psi^{-1}=\psi\circ\sigma\circ
\psi^{-1},
\end{align*}
this rewrites as $\left(  -1\right)  _{\psi\circ f^{-1}}^{f\circ\sigma\circ
f^{-1}}=\left(  -1\right)  ^{\psi\circ\sigma\circ\psi^{-1}}$.
\end{vershort}

\begin{verlong}
On the other hand, the map $\psi\circ f^{-1}:Y\rightarrow\left\{
1,2,\ldots,n\right\}  $ is bijective (since it is the composition of the two
bijective maps $\psi$ and $f^{-1}$). In other words, $\psi\circ f^{-1}$ is a
bijection from $Y$ to $\left\{  1,2,\ldots,n\right\}  $. Hence, the definition
of $\left(  -1\right)  _{\psi\circ f^{-1}}^{f\circ\sigma\circ f^{-1}}$ (given
in Exercise \ref{exe.ps4.2}) shows that
\[
\left(  -1\right)  _{\psi\circ f^{-1}}^{f\circ\sigma\circ f^{-1}}=\left(
-1\right)  ^{\left(  \psi\circ f^{-1}\right)  \circ\left(  f\circ\sigma\circ
f^{-1}\right)  \circ\left(  \psi\circ f^{-1}\right)  ^{-1}}.
\]
In view of
\begin{align*}
&  \underbrace{\left(  \psi\circ f^{-1}\right)  \circ\left(  f\circ\sigma\circ
f^{-1}\right)  }_{=\psi\circ f^{-1}\circ f\circ\sigma\circ f^{-1}}%
\circ\underbrace{\left(  \psi\circ f^{-1}\right)  ^{-1}}_{=\left(
f^{-1}\right)  ^{-1}\circ\psi^{-1}}\\
&  =\psi\circ\underbrace{f^{-1}\circ f}_{=\operatorname*{id}\nolimits_{X}%
}\circ\sigma\circ\underbrace{f^{-1}\circ\left(  f^{-1}\right)  ^{-1}%
}_{=\operatorname*{id}\nolimits_{X}}\circ\psi^{-1}=\psi\circ\sigma\circ
\psi^{-1},
\end{align*}
this rewrites as $\left(  -1\right)  _{\psi\circ f^{-1}}^{f\circ\sigma\circ
f^{-1}}=\left(  -1\right)  ^{\psi\circ\sigma\circ\psi^{-1}}$.
\end{verlong}

But the definition of $\left(  -1\right)  ^{\sigma}$ (given in Exercise
\ref{exe.ps4.2}) shows that $\left(  -1\right)  ^{\sigma}=\left(  -1\right)
_{\phi}^{\sigma}$ for any bijection $\phi:X\rightarrow\left\{  1,2,\ldots
,n\right\}  $. Applying this to $\phi=\psi$, we obtain
\begin{equation}
\left(  -1\right)  ^{\sigma}=\left(  -1\right)  _{\psi}^{\sigma}=\left(
-1\right)  ^{\psi\circ\sigma\circ\psi^{-1}}.
\label{pf.prop.sol.sigmacrosstau.signs-conj.1}%
\end{equation}

Also, the definition of $\left(  -1\right)  ^{f\circ\sigma\circ f^{-1}}$
(given in Exercise \ref{exe.ps4.2}) shows that $\left(  -1\right)
^{f\circ\sigma\circ f^{-1}}=\left(  -1\right)  _{\phi}^{f\circ\sigma\circ
f^{-1}}$ for any bijection $\phi:Y\rightarrow\left\{  1,2,\ldots,n\right\}  $.
Applying this to $\phi=\psi\circ f^{-1}$, we obtain
\[
\left(  -1\right)  ^{f\circ\sigma\circ f^{-1}}=\left(  -1\right)  _{\psi\circ
f^{-1}}^{f\circ\sigma\circ f^{-1}}=\left(  -1\right)  ^{\psi\circ\sigma
\circ\psi^{-1}}.
\]
Comparing this with (\ref{pf.prop.sol.sigmacrosstau.signs-conj.1}), we obtain
$\left(  -1\right)  ^{f\circ\sigma\circ f^{-1}}=\left(  -1\right)  ^{\sigma}$.
This completes the proof of Proposition
\ref{prop.sol.sigmacrosstau.signs-conj}.
\end{proof}

\begin{corollary}
\label{cor.sol.sigmacrosstau.signs-fggf}Let $U$ and $V$ be two finite sets.
Let $\sigma$ be a permutation of $U$. Let $\tau$ be a permutation of $V$.
Then, $\left(  -1\right)  ^{\sigma\times\tau}=\left(  -1\right)  ^{\tau
\times\sigma}$.
\end{corollary}

\begin{proof}
[Proof of Corollary \ref{cor.sol.sigmacrosstau.signs-fggf}.]Corollary
\ref{cor.sol.sigmacrosstau.permxperm} yields that $\sigma\times\tau$ is a
permutation of $U\times V$. Hence, $\left(  -1\right)  ^{\sigma\times\tau}$ is well-defined.

\begin{verlong}
Also, Corollary \ref{cor.sol.sigmacrosstau.permxperm} (applied to $V$, $U$,
$\tau$ and $\sigma$ instead of $U$, $V$, $\sigma$ and $\tau$) shows that
$\tau\times\sigma$ is a permutation of $V\times U$. Hence, $\left(  -1\right)
^{\tau\times\sigma}$ is well-defined.
\end{verlong}

Consider the map $\gamma:U\times V\rightarrow V\times U$ defined in
Proposition \ref{prop.sol.sigmacrosstau.trans}. Then, Proposition
\ref{prop.sol.sigmacrosstau.trans} \textbf{(a)} shows that this map $\gamma$
is a bijection. Hence, Proposition \ref{prop.sol.sigmacrosstau.signs-conj}
(applied to $U\times V$, $V\times U$, $\gamma$ and $\sigma\times\tau$ instead
of $X$, $Y$, $f$ and $\sigma$) shows that $\gamma\circ\left(  \sigma\times
\tau\right)  \circ\gamma^{-1}$ is a permutation of $V\times U$ and satisfies
$\left(  -1\right)  ^{\gamma\circ\left(  \sigma\times\tau\right)  \circ
\gamma^{-1}}=\left(  -1\right)  ^{\sigma\times\tau}$.

But Proposition \ref{prop.sol.sigmacrosstau.trans} \textbf{(b)} (applied to
$f=\sigma$ and $g=\tau$) shows that $\sigma\times\tau=\gamma^{-1}\circ\left(
\tau\times\sigma\right)  \circ\gamma$. Hence,%
\[
\gamma\circ\underbrace{\left(  \sigma\times\tau\right)  }_{=\gamma^{-1}%
\circ\left(  \tau\times\sigma\right)  \circ\gamma}\circ\gamma^{-1}%
=\underbrace{\gamma\circ\gamma^{-1}}_{=\operatorname*{id}\nolimits_{V\times
U}}\circ\left(  \tau\times\sigma\right)  \circ\underbrace{\gamma\circ
\gamma^{-1}}_{=\operatorname*{id}\nolimits_{V\times U}}=\tau\times\sigma.
\]
Thus, $\left(  -1\right)  ^{\gamma\circ\left(  \sigma\times\tau\right)
\circ\gamma^{-1}}=\left(  -1\right)  ^{\tau\times\sigma}$. Comparing this with
$\left(  -1\right)  ^{\gamma\circ\left(  \sigma\times\tau\right)  \circ
\gamma^{-1}}=\left(  -1\right)  ^{\sigma\times\tau}$, we obtain $\left(
-1\right)  ^{\sigma\times\tau}=\left(  -1\right)  ^{\tau\times\sigma}$. This
proves Corollary \ref{cor.sol.sigmacrosstau.signs-fggf}.
\end{proof}

For the rest of Section \ref{sect.sol.perm.sigmacrosstau}, we shall use
Definition \ref{def.transposX}.

\begin{proposition}
\label{prop.sol.sigmacrosstau.transposX.cross1}Let $X$ be a set. Let
$m\in\mathbb{N}$. Let $M$ be the set $\left\{  1,2,\ldots,m\right\}  $. Let
$x$ and $y$ be two distinct elements of $X$. Then, $t_{\left(  x,1\right)
,\left(  y,1\right)  },t_{\left(  x,2\right)  ,\left(  y,2\right)  }%
,\ldots,t_{\left(  x,m\right)  ,\left(  y,m\right)  }$ are $m$ well-defined
transpositions of the set $X\times M$, and we have%
\[
t_{x,y}\times\operatorname*{id}\nolimits_{M}=t_{\left(  x,m\right)  ,\left(
y,m\right)  }\circ t_{\left(  x,m-1\right)  ,\left(  y,m-1\right)  }%
\circ\cdots\circ t_{\left(  x,1\right)  ,\left(  y,1\right)  }.
\]

\end{proposition}

\begin{vershort}
\begin{proof}
[Proof of Proposition \ref{prop.sol.sigmacrosstau.transposX.cross1}.]Let us
first give a quick informal proof of Proposition
\ref{prop.sol.sigmacrosstau.transposX.cross1}:

For each $j\in\left\{  1,2,\ldots,m\right\}  $, the transposition $t_{\left(
x,j\right)  ,\left(  y,j\right)  }$ swaps $\left(  x,j\right)  $ with $\left(
y,j\right)  $ while leaving all other elements of $X\times M$ unchanged. Thus,
the composition $t_{\left(  x,m\right)  ,\left(  y,m\right)  }\circ t_{\left(
x,m-1\right)  ,\left(  y,m-1\right)  }\circ\cdots\circ t_{\left(  x,1\right)
,\left(  y,1\right)  }$ of these $m$ transpositions \newline$t_{\left(
x,m\right)  ,\left(  y,m\right)  },t_{\left(  x,m-1\right)  ,\left(
y,m-1\right)  },\ldots,t_{\left(  x,1\right)  ,\left(  y,1\right)  }$ swaps
the elements $\left(  x,1\right)  ,\left(  x,2\right)  ,\ldots,\left(
x,m\right)  $ with the elements $\left(  y,1\right)  ,\left(  y,2\right)
,\ldots,\left(  y,m\right)  $, respectively, while leaving all other elements
of $X\times M$ unchanged\footnote{Here we are using the observation that the
$2m$ elements $\left(  x,1\right)  ,\left(  y,1\right)  ,\left(  x,2\right)
,\left(  y,2\right)  ,\ldots,\left(  x,m\right)  ,\left(  y,m\right)  $ of
$X\times M$ are distinct (since $x$ and $y$ are distinct), whence any element
of $X\times M$ is moved (i.e., not left unchanged) by \textbf{at most one} of
the $m$ transpositions $t_{\left(  x,m\right)  ,\left(  y,m\right)
},t_{\left(  x,m-1\right)  ,\left(  y,m-1\right)  },\ldots,t_{\left(
x,1\right)  ,\left(  y,1\right)  }$.}. In other words, $t_{\left(  x,m\right)
,\left(  y,m\right)  }\circ t_{\left(  x,m-1\right)  ,\left(  y,m-1\right)
}\circ\cdots\circ t_{\left(  x,1\right)  ,\left(  y,1\right)  }$ swaps each
element of the form $\left(  x,q\right)  $ (with $q\in\left\{  1,2,\ldots
,m\right\}  $) with the corresponding $\left(  y,q\right)  $, while leaving
all other elements of $X\times M$ unchanged. In other words,%
\begin{align}
&  \left(  t_{\left(  x,m\right)  ,\left(  y,m\right)  }\circ t_{\left(
x,m-1\right)  ,\left(  y,m-1\right)  }\circ\cdots\circ t_{\left(  x,1\right)
,\left(  y,1\right)  }\right)  \left(  \left(  u,q\right)  \right) \nonumber\\
&  =%
\begin{cases}
\left(  y,q\right)  , & \text{if }u=x;\\
\left(  x,q\right)  , & \text{if }u=y;\\
\left(  u,q\right)  , & \text{otherwise}%
\end{cases}
\label{pf.prop.sol.sigmacrosstau.transposX.cross1.short.1}%
\end{align}
for each $\left(  u,q\right)  \in X\times M$.

But the transposition $t_{x,y}$ swaps $x$ with $y$ while leaving all other
elements of $X$ unchanged. Thus,
\begin{equation}
t_{x,y}\left(  u\right)  =%
\begin{cases}
y, & \text{if }u=x;\\
x, & \text{if }u=y;\\
u, & \text{otherwise}%
\end{cases}
\label{pf.prop.sol.sigmacrosstau.transposX.cross1.short.2}%
\end{equation}
for each $u\in X$. Hence, for each $\left(  u,q\right)  \in X\times M$, we
have%
\begin{align*}
\left(  t_{x,y}\times\operatorname*{id}\nolimits_{M}\right)  \left(  \left(
u,q\right)  \right)   &  =\left(  t_{x,y}\left(  u\right)
,\underbrace{\operatorname*{id}\nolimits_{M}\left(  q\right)  }_{=q}\right)
=\left(  t_{x,y}\left(  u\right)  ,q\right) \\
&  =\left(
\begin{cases}
y, & \text{if }u=x;\\
x, & \text{if }u=y;\\
u, & \text{otherwise}%
\end{cases}
,q\right)  \ \ \ \ \ \ \ \ \ \ \left(  \text{by
(\ref{pf.prop.sol.sigmacrosstau.transposX.cross1.short.2})}\right) \\
&  =%
\begin{cases}
\left(  y,q\right)  , & \text{if }u=x;\\
\left(  x,q\right)  , & \text{if }u=y;\\
\left(  u,q\right)  , & \text{otherwise}%
\end{cases}
\\
&  =\left(  t_{\left(  x,m\right)  ,\left(  y,m\right)  }\circ t_{\left(
x,m-1\right)  ,\left(  y,m-1\right)  }\circ\cdots\circ t_{\left(  x,1\right)
,\left(  y,1\right)  }\right)  \left(  \left(  u,q\right)  \right)
\end{align*}
(by (\ref{pf.prop.sol.sigmacrosstau.transposX.cross1.short.1})). In other
words, we have%
\[
t_{x,y}\times\operatorname*{id}\nolimits_{M}=t_{\left(  x,m\right)  ,\left(
y,m\right)  }\circ t_{\left(  x,m-1\right)  ,\left(  y,m-1\right)  }%
\circ\cdots\circ t_{\left(  x,1\right)  ,\left(  y,1\right)  }.
\]
Thus, Proposition \ref{prop.sol.sigmacrosstau.transposX.cross1} follows.

Let us now show how to prove Proposition
\ref{prop.sol.sigmacrosstau.transposX.cross1} rigorously. The following
rigorous proof is, of course, just a formalization of the above informal argument.

We have $x\neq y$ (since $x$ and $y$ are distinct). Thus, for each
$p\in\left\{  1,2,\ldots,m\right\}  $, we have $\left(  x,p\right)
\neq\left(  y,p\right)  $. Hence, for each $p\in\left\{  1,2,\ldots,m\right\}
$, the transposition $t_{\left(  x,p\right)  ,\left(  y,p\right)  }$ of
$X\times M$ is well-defined. In other words, $t_{\left(  x,1\right)  ,\left(
y,1\right)  },t_{\left(  x,2\right)  ,\left(  y,2\right)  },\ldots,t_{\left(
x,m\right)  ,\left(  y,m\right)  }$ are $m$ well-defined transpositions of the
set $X\times M$.

It remains to prove that $t_{x,y}\times\operatorname*{id}\nolimits_{M}%
=t_{\left(  x,m\right)  ,\left(  y,m\right)  }\circ t_{\left(  x,m-1\right)
,\left(  y,m-1\right)  }\circ\cdots\circ t_{\left(  x,1\right)  ,\left(
y,1\right)  }$.

Let us first observe the following fact: For any $j\in\left\{  1,2,\ldots
,m\right\}  $ and $q\in M$ and $z\in X$ satisfying $j\neq q$, we have%
\begin{equation}
t_{\left(  x,j\right)  ,\left(  y,j\right)  }\left(  \left(  z,q\right)
\right)  =\left(  z,q\right)  .
\label{pf.prop.sol.sigmacrosstau.transposX.cross1.short.t-dist}%
\end{equation}

[\textit{Proof of
(\ref{pf.prop.sol.sigmacrosstau.transposX.cross1.short.t-dist}):} Let
$j\in\left\{  1,2,\ldots,m\right\}  $ and $q\in M$ and $z\in X$ be such that
$j\neq q$.

The pair $\left(  z,q\right)  $ is distinct from both $\left(  x,j\right)  $
and $\left(  y,j\right)  $ (since $j\neq q$). In other words, $\left(
z,q\right)  \notin\left\{  \left(  x,j\right)  ,\left(  y,j\right)  \right\}
$. Combining this with $\left(  z,q\right)  \in X\times M$ (since $z\in X$ and
$q\in M$), we obtain $\left(  z,q\right)  \in\left(  X\times M\right)
\setminus\left\{  \left(  x,j\right)  ,\left(  y,j\right)  \right\}  $. Hence,
Lemma \ref{lem.sol.perm.transX.tij1} \textbf{(c)} (applied to $X\times M$,
$\left(  x,j\right)  $, $\left(  y,j\right)  $ and $\left(  z,q\right)  $
instead of $X$, $i$, $j$ and $k$) yields $t_{\left(  x,j\right)  ,\left(
y,j\right)  }\left(  \left(  z,q\right)  \right)  =\left(  z,q\right)  $. This
proves (\ref{pf.prop.sol.sigmacrosstau.transposX.cross1.short.t-dist}).]

Define a map $s:X\times M\rightarrow X\times M$ by%
\begin{equation}
s=t_{\left(  x,m\right)  ,\left(  y,m\right)  }\circ t_{\left(  x,m-1\right)
,\left(  y,m-1\right)  }\circ\cdots\circ t_{\left(  x,1\right)  ,\left(
y,1\right)  }. \label{pf.prop.sol.sigmacrosstau.transposX.cross1.short.s=}%
\end{equation}

Fix $c\in X\times M$. We shall show that $\left(  t_{x,y}\times
\operatorname*{id}\nolimits_{M}\right)  \left(  c\right)  =s\left(  c\right)
$.

Write the element $c\in X\times M$ in the form $c=\left(  z,q\right)  $ for
some $z\in X$ and $q\in M$.

Applying the map $t_{x,y}\times\operatorname*{id}\nolimits_{M}$ to both sides
of the equality $c=\left(  z,q\right)  $, we obtain%
\begin{align}
\left(  t_{x,y}\times\operatorname*{id}\nolimits_{M}\right)  \left(  c\right)
&  =\left(  t_{x,y}\times\operatorname*{id}\nolimits_{M}\right)  \left(
\left(  z,q\right)  \right)  =\left(  t_{x,y}\left(  z\right)
,\underbrace{\operatorname*{id}\nolimits_{M}\left(  q\right)  }_{=q}\right)
\nonumber\\
&  =\left(  t_{x,y}\left(  z\right)  ,q\right)  .
\label{pf.prop.sol.sigmacrosstau.transposX.cross1.short.LHS1}%
\end{align}

We have $q\in M=\left\{  1,2,\ldots,m\right\}  $. We are in one of the
following three cases:

\textit{Case 1:} We have $z=x$.

\textit{Case 2:} We have $z=y$.

\textit{Case 3:} We have neither $z=x$ nor $z=y$.

Let us first consider Case 1. In this case, we have $z=x$. But Lemma
\ref{lem.sol.perm.transX.tij1} \textbf{(a)} (applied to $x$ and $y$ instead of
$i$ and $j$) yields $t_{x,y}\left(  x\right)  =y$. But $t_{x,y}\left(
\underbrace{z}_{=x}\right)  =t_{x,y}\left(  x\right)  =y$.

Hence, (\ref{pf.prop.sol.sigmacrosstau.transposX.cross1.short.LHS1}) yields
\begin{equation}
\left(  t_{x,y}\times\operatorname*{id}\nolimits_{M}\right)  \left(  c\right)
=\left(  \underbrace{t_{x,y}\left(  z\right)  }_{=y},q\right)  =\left(
y,q\right)  . \label{pf.prop.sol.sigmacrosstau.transposX.cross1.short.c1.LHS}%
\end{equation}

On the other hand, $\left(  x,q\right)  \neq\left(  y,q\right)  $ (since
$x\neq y$). Hence, the two elements $\left(  x,q\right)  $ and $\left(
y,q\right)  $ of $X\times M$ are distinct. Thus, Lemma
\ref{lem.sol.perm.transX.tij1} \textbf{(a)} (applied to $X\times M$, $\left(
x,q\right)  $ and $\left(  y,q\right)  $ instead of $X$, $i$ and $j$) yields
\begin{equation}
t_{\left(  x,q\right)  ,\left(  y,q\right)  }\left(  \left(  x,q\right)
\right)  =\left(  y,q\right)  .
\label{pf.prop.sol.sigmacrosstau.transposX.cross1.short.c1.0}%
\end{equation}

Next, we observe that%
\begin{equation}
t_{\left(  x,j\right)  ,\left(  y,j\right)  }\left(  \left(  x,q\right)
\right)  =\left(  x,q\right)  \ \ \ \ \ \ \ \ \ \ \text{for each }j\in\left\{
1,2,\ldots,m\right\}  \text{ satisfying }j<q.
\label{pf.prop.sol.sigmacrosstau.transposX.cross1.short.c1.1}%
\end{equation}

[\textit{Proof of (\ref{pf.prop.sol.sigmacrosstau.transposX.cross1.short.c1.1}%
):} Let $j\in\left\{  1,2,\ldots,m\right\}  $ be such that $j<q$. Thus, $j\neq
q$. Hence, (\ref{pf.prop.sol.sigmacrosstau.transposX.cross1.short.t-dist})
(applied to $z=x$) yields $t_{\left(  x,j\right)  ,\left(  y,j\right)
}\left(  \left(  x,q\right)  \right)  =\left(  x,q\right)  $. This proves
(\ref{pf.prop.sol.sigmacrosstau.transposX.cross1.short.c1.1}).]

Next, we observe that%
\begin{equation}
t_{\left(  x,j\right)  ,\left(  y,j\right)  }\left(  \left(  y,q\right)
\right)  =\left(  y,q\right)  \ \ \ \ \ \ \ \ \ \ \text{for each }j\in\left\{
1,2,\ldots,m\right\}  \text{ satisfying }j>q.
\label{pf.prop.sol.sigmacrosstau.transposX.cross1.short.c1.2}%
\end{equation}

[\textit{Proof of (\ref{pf.prop.sol.sigmacrosstau.transposX.cross1.short.c1.2}%
):} Let $j\in\left\{  1,2,\ldots,m\right\}  $ be such that $j>q$. Thus, $j\neq
q$. Hence, (\ref{pf.prop.sol.sigmacrosstau.transposX.cross1.short.t-dist})
(applied to $z=y$) yields $t_{\left(  x,j\right)  ,\left(  y,j\right)
}\left(  \left(  y,q\right)  \right)  =\left(  y,q\right)  $. This proves
(\ref{pf.prop.sol.sigmacrosstau.transposX.cross1.short.c1.2}).]

Now, we can apply Lemma \ref{lem.sol.sigmacrosstau.compose1} to $X\times M$,
$t_{\left(  x,j\right)  ,\left(  y,j\right)  }$, $\left(  x,q\right)  $,
$\left(  y,q\right)  $ and $q$ instead of $X$, $f_{j}$, $x$, $y$ and $i$
(since the equalities
(\ref{pf.prop.sol.sigmacrosstau.transposX.cross1.short.c1.0}),
(\ref{pf.prop.sol.sigmacrosstau.transposX.cross1.short.c1.1}) and
(\ref{pf.prop.sol.sigmacrosstau.transposX.cross1.short.c1.2}) hold). Thus, we
obtain%
\begin{equation}
\left(  t_{\left(  x,m\right)  ,\left(  y,m\right)  }\circ t_{\left(
x,m-1\right)  ,\left(  y,m-1\right)  }\circ\cdots\circ t_{\left(  x,1\right)
,\left(  y,1\right)  }\right)  \left(  \left(  x,q\right)  \right)  =\left(
y,q\right)  . \label{pf.prop.sol.sigmacrosstau.transposX.cross1.short.c1.5}%
\end{equation}
In view of (\ref{pf.prop.sol.sigmacrosstau.transposX.cross1.short.s=}), this
rewrites as $s\left(  \left(  x,q\right)  \right)  =\left(  y,q\right)  $.
Comparing this with
(\ref{pf.prop.sol.sigmacrosstau.transposX.cross1.short.c1.LHS}), we obtain
$\left(  t_{x,y}\times\operatorname*{id}\nolimits_{M}\right)  \left(
c\right)  =s\left(  \left(  x,q\right)  \right)  $. Comparing this with
$s\left(  \underbrace{c}_{=\left(  z,q\right)  }\right)  =s\left(  \left(
\underbrace{z}_{=x},q\right)  \right)  =s\left(  \left(  x,q\right)  \right)
$, we obtain $\left(  t_{x,y}\times\operatorname*{id}\nolimits_{M}\right)
\left(  c\right)  =s\left(  c\right)  $. Hence, $\left(  t_{x,y}%
\times\operatorname*{id}\nolimits_{M}\right)  \left(  c\right)  =s\left(
c\right)  $ is proven in Case 1.

We leave it to the reader to verify $\left(  t_{x,y}\times\operatorname*{id}%
\nolimits_{M}\right)  \left(  c\right)  =s\left(  c\right)  $ in Case 2. (The
verification is analogous to what we did in Case 1, but the roles of $x$ and
$y$ are often interchanged, and instead of Lemma
\ref{lem.sol.perm.transX.tij1} \textbf{(a)} we must now apply Lemma
\ref{lem.sol.perm.transX.tij1} \textbf{(b)}.)

Let us finally consider Case 3. In this case, we have neither $z=x$ nor $z=y$.
Thus, $z\notin\left\{  x,y\right\}  $. Combining this with $z\in X$, we obtain
$z\in X\setminus\left\{  x,y\right\}  $. Thus, Lemma
\ref{lem.sol.perm.transX.tij1} \textbf{(c)} (applied to $x$, $y$ and $z$
instead of $i$, $j$ and $k$) yields $t_{x,y}\left(  z\right)  =z$.

Hence, (\ref{pf.prop.sol.sigmacrosstau.transposX.cross1.short.LHS1}) yields
\begin{equation}
\left(  t_{x,y}\times\operatorname*{id}\nolimits_{M}\right)  \left(  c\right)
=\left(  \underbrace{t_{x,y}\left(  z\right)  }_{=z},q\right)  =\left(
z,q\right)  . \label{pf.prop.sol.sigmacrosstau.transposX.cross1.short.c3.LHS}%
\end{equation}

Next, we observe that%
\begin{equation}
t_{\left(  x,j\right)  ,\left(  y,j\right)  }\left(  \left(  z,q\right)
\right)  =\left(  z,q\right)  \ \ \ \ \ \ \ \ \ \ \text{for each }j\in\left\{
1,2,\ldots,m\right\}  .
\label{pf.prop.sol.sigmacrosstau.transposX.cross1.short.c3.1}%
\end{equation}

[\textit{Proof of (\ref{pf.prop.sol.sigmacrosstau.transposX.cross1.short.c3.1}%
):} Let $j\in\left\{  1,2,\ldots,m\right\}  $. Thus, $j\in\left\{
1,2,\ldots,m\right\}  =M$. The elements $\left(  x,j\right)  $ and $\left(
y,j\right)  $ of $X\times M$ are distinct (since $x\neq y$). The element
$\left(  z,q\right)  $ is distinct from both $\left(  x,j\right)  $ and
$\left(  y,j\right)  $ (since we have neither $z=x$ nor $z=y$). In other
words, $\left(  z,q\right)  \notin\left\{  \left(  x,j\right)  ,\left(
y,j\right)  \right\}  $. Combining this with $\left(  z,q\right)  \in X\times
M$, we obtain $\left(  z,q\right)  \in\left(  X\times M\right)  \setminus
\left\{  \left(  x,j\right)  ,\left(  y,j\right)  \right\}  $. Hence, Lemma
\ref{lem.sol.perm.transX.tij1} \textbf{(c)} (applied to $X\times M$, $\left(
x,j\right)  $, $\left(  y,j\right)  $ and $\left(  z,q\right)  $ instead of
$X$, $i$, $j$ and $k$) yields $t_{\left(  x,j\right)  ,\left(  y,j\right)
}\left(  \left(  z,q\right)  \right)  =\left(  z,q\right)  $. This proves
(\ref{pf.prop.sol.sigmacrosstau.transposX.cross1.short.c3.1}).]

Now, we can apply Lemma \ref{lem.sol.sigmacrosstau.compose0} to $X\times M$,
$t_{\left(  x,j\right)  ,\left(  y,j\right)  }$ and $\left(  z,q\right)  $
instead of $X$, $f_{j}$ and $x$ (since
(\ref{pf.prop.sol.sigmacrosstau.transposX.cross1.short.c3.1}) holds). Thus, we
obtain%
\begin{equation}
\left(  t_{\left(  x,m\right)  ,\left(  y,m\right)  }\circ t_{\left(
x,m-1\right)  ,\left(  y,m-1\right)  }\circ\cdots\circ t_{\left(  x,1\right)
,\left(  y,1\right)  }\right)  \left(  \left(  z,q\right)  \right)  =\left(
z,q\right)  . \label{pf.prop.sol.sigmacrosstau.transposX.cross1.short.c3.5}%
\end{equation}
In view of (\ref{pf.prop.sol.sigmacrosstau.transposX.cross1.short.s=}), this
rewrites as $s\left(  \left(  z,q\right)  \right)  =\left(  z,q\right)  $.
Comparing this with
(\ref{pf.prop.sol.sigmacrosstau.transposX.cross1.short.c3.LHS}), we obtain
$\left(  t_{x,y}\times\operatorname*{id}\nolimits_{M}\right)  \left(
c\right)  =s\left(  \left(  z,q\right)  \right)  $. Comparing this with
$s\left(  \underbrace{c}_{=\left(  z,q\right)  }\right)  =s\left(  \left(
z,q\right)  \right)  $, we obtain $\left(  t_{x,y}\times\operatorname*{id}%
\nolimits_{M}\right)  \left(  c\right)  =s\left(  c\right)  $. Hence, $\left(
t_{x,y}\times\operatorname*{id}\nolimits_{M}\right)  \left(  c\right)
=s\left(  c\right)  $ is proven in Case 3.

We have now proven $\left(  t_{x,y}\times\operatorname*{id}\nolimits_{M}%
\right)  \left(  c\right)  =s\left(  c\right)  $ in each of the three Cases 1,
2 and 3. Thus, $\left(  t_{x,y}\times\operatorname*{id}\nolimits_{M}\right)
\left(  c\right)  =s\left(  c\right)  $ always holds.

Now, forget that we fixed $c$. We thus have proven that $\left(  t_{x,y}%
\times\operatorname*{id}\nolimits_{M}\right)  \left(  c\right)  =s\left(
c\right)  $ for each $c\in X\times M$. Thus,%
\[
t_{x,y}\times\operatorname*{id}\nolimits_{M}=s=t_{\left(  x,m\right)  ,\left(
y,m\right)  }\circ t_{\left(  x,m-1\right)  ,\left(  y,m-1\right)  }%
\circ\cdots\circ t_{\left(  x,1\right)  ,\left(  y,1\right)  }.
\]
This completes the rigorous proof of Proposition
\ref{prop.sol.sigmacrosstau.transposX.cross1}.
\end{proof}
\end{vershort}

\begin{verlong}
\begin{proof}
[Proof of Proposition \ref{prop.sol.sigmacrosstau.transposX.cross1}.]For each
$p\in\left\{  1,2,\ldots,m\right\}  $, the transposition $t_{\left(
x,p\right)  ,\left(  y,p\right)  }$ of $X\times M$ is
well-defined\footnote{\textit{Proof.} Let $p\in\left\{  1,2,\ldots,m\right\}
$. We must prove that the transposition $t_{\left(  x,p\right)  ,\left(
y,p\right)  }$ of $X\times M$ is well-defined.
\par
Assume (for the sake of contradiction) that $\left(  x,p\right)  =\left(
y,p\right)  $. Thus, $x=y$ and $p=p$. But $x=y$ contradicts the fact that $x$
and $y$ are distinct. This contradiction shows that our assumption (that
$\left(  x,p\right)  =\left(  y,p\right)  $) was wrong. Hence, we cannot have
$\left(  x,p\right)  =\left(  y,p\right)  $. In other words, we have $\left(
x,p\right)  \neq\left(  y,p\right)  $.
\par
We have $\left(  x,p\right)  \in X\times M$ (since $x\in X$ and $p\in\left\{
1,2,\ldots,m\right\}  =M$). The same argument (applied to $y$ instead of $x$)
yields $\left(  y,p\right)  \in X\times M$. Hence, $\left(  x,p\right)  $ and
$\left(  y,p\right)  $ are two elements of $X\times M$. These two elements are
distinct (since $\left(  x,p\right)  \neq\left(  y,p\right)  $). Thus, the
transposition $t_{\left(  x,p\right)  ,\left(  y,p\right)  }$ of $X\times M$
is well-defined (according to Definition \ref{def.transposX}). Qed.}. In other
words, $t_{\left(  x,1\right)  ,\left(  y,1\right)  },t_{\left(  x,2\right)
,\left(  y,2\right)  },\ldots,t_{\left(  x,m\right)  ,\left(  y,m\right)  }$
are $m$ well-defined transpositions of the set $X\times M$. Thus, $t_{\left(
x,1\right)  ,\left(  y,1\right)  },t_{\left(  x,2\right)  ,\left(  y,2\right)
},\ldots,t_{\left(  x,m\right)  ,\left(  y,m\right)  }$ are $m$ maps $X\times
M\rightarrow X\times Y$.

It remains to prove that $t_{x,y}\times\operatorname*{id}\nolimits_{M}%
=t_{\left(  x,m\right)  ,\left(  y,m\right)  }\circ t_{\left(  x,m-1\right)
,\left(  y,m-1\right)  }\circ\cdots\circ t_{\left(  x,1\right)  ,\left(
y,1\right)  }$.

Let us first observe the following fact: For any $j\in\left\{  1,2,\ldots
,m\right\}  $ and $q\in M$ and $z\in X$ satisfying $j\neq q$, we have%
\begin{equation}
t_{\left(  x,j\right)  ,\left(  y,j\right)  }\left(  \left(  z,q\right)
\right)  =\left(  z,q\right)  .
\label{pf.prop.sol.sigmacrosstau.transposX.cross1.t-dist}%
\end{equation}

[\textit{Proof of (\ref{pf.prop.sol.sigmacrosstau.transposX.cross1.t-dist}):}
Let $j\in\left\{  1,2,\ldots,m\right\}  $ and $q\in M$ and $z\in X$ be such
that $j\neq q$.

We have $\left(  z,q\right)  \neq\left(  w,j\right)  $ for each $w\in
X$\ \ \ \footnote{\textit{Proof.} Let $w\in X$. We must show that $\left(
z,q\right)  \neq\left(  w,j\right)  $.
\par
Indeed, assume the contrary. Thus, $\left(  z,q\right)  =\left(  w,j\right)
$. In other words, $\left(  w,j\right)  =\left(  z,q\right)  $. In other
words, $w=z$ and $j=q$. But $j=q$ contradicts $j\neq q$. This contradiction
shows that our assumption was false.
\par
Hence, $\left(  z,q\right)  \neq\left(  w,j\right)  $ is proven.}. Applying
this to $w=x$, we obtain $\left(  z,q\right)  \neq\left(  x,j\right)  $. In
other words, we don't have $\left(  z,q\right)  =\left(  x,j\right)  $.

Recall that $\left(  z,q\right)  \neq\left(  w,j\right)  $ for each $w\in X$.
Applying this to $w=y$, we obtain $\left(  z,q\right)  \neq\left(  y,j\right)
$. In other words, we don't have $\left(  z,q\right)  =\left(  y,j\right)  $.
Hence, we have neither $\left(  z,q\right)  =\left(  x,j\right)  $ nor
$\left(  z,q\right)  =\left(  y,j\right)  $.

We have $j\in\left\{  1,2,\ldots,m\right\}  =M$. The elements $\left(
x,j\right)  $ and $\left(  y,j\right)  $ of $X\times M$ are
distinct\footnote{\textit{Proof.} Assume the contrary. Thus, $\left(
x,j\right)  =\left(  y,j\right)  $. In other words, $x=y$ and $j=j$. But $x=y$
contradicts the fact that $x$ and $y$ are distinct. This contradiction shows
that our assumption was wrong. Qed.}.

We have $\left(  z,q\right)  \notin\left\{  \left(  x,j\right)  ,\left(
y,j\right)  \right\}  $\ \ \ \ \footnote{\textit{Proof.} Assume the contrary.
Thus, $\left(  z,q\right)  \in\left\{  \left(  x,j\right)  ,\left(
y,j\right)  \right\}  $. In other words, $\left(  z,q\right)  $ equals either
$\left(  x,j\right)  $ or $\left(  y,j\right)  $. In other words, we have
either $\left(  z,q\right)  =\left(  x,j\right)  $ or $\left(  z,q\right)
=\left(  y,j\right)  $. This contradicts the fact that we have neither
$\left(  z,q\right)  =\left(  x,j\right)  $ nor $\left(  z,q\right)  =\left(
y,j\right)  $. This contradiction shows that our assumption was wrong. Qed.}.
Combining $\left(  z,q\right)  \in X\times M$ (since $z\in X$ and $q\in M$)
with $\left(  z,q\right)  \notin\left\{  \left(  x,j\right)  ,\left(
y,j\right)  \right\}  $, we obtain $\left(  z,q\right)  \in\left(  X\times
M\right)  \setminus\left\{  \left(  x,j\right)  ,\left(  y,j\right)  \right\}
$. Hence, Lemma \ref{lem.sol.perm.transX.tij1} \textbf{(c)} (applied to
$X\times M$, $\left(  x,j\right)  $, $\left(  y,j\right)  $ and $\left(
z,q\right)  $ instead of $X$, $i$, $j$ and $k$) yields $t_{\left(  x,j\right)
,\left(  y,j\right)  }\left(  \left(  z,q\right)  \right)  =\left(
z,q\right)  $. This proves
(\ref{pf.prop.sol.sigmacrosstau.transposX.cross1.t-dist}).]

Define a map $s:X\times M\rightarrow X\times M$ by%
\begin{equation}
s=t_{\left(  x,m\right)  ,\left(  y,m\right)  }\circ t_{\left(  x,m-1\right)
,\left(  y,m-1\right)  }\circ\cdots\circ t_{\left(  x,1\right)  ,\left(
y,1\right)  }. \label{pf.prop.sol.sigmacrosstau.transposX.cross1.s=}%
\end{equation}

Fix $c\in X\times M$. We shall show that $\left(  t_{x,y}\times
\operatorname*{id}\nolimits_{M}\right)  \left(  c\right)  =s\left(  c\right)
$.

We have $c\in X\times M$. In other words, we can write $c$ in the form
$c=\left(  z,q\right)  $ for some $z\in X$ and $q\in M$. Consider these $z$
and $q$.

Applying the map $t_{x,y}\times\operatorname*{id}\nolimits_{M}$ to both sides
of the equality $c=\left(  z,q\right)  $, we obtain%
\begin{align}
\left(  t_{x,y}\times\operatorname*{id}\nolimits_{M}\right)  \left(  c\right)
&  =\left(  t_{x,y}\times\operatorname*{id}\nolimits_{M}\right)  \left(
\left(  z,q\right)  \right) \nonumber\\
&  =\left(  t_{x,y}\left(  z\right)  ,\underbrace{\operatorname*{id}%
\nolimits_{M}\left(  q\right)  }_{=q}\right)  \ \ \ \ \ \ \ \ \ \ \left(
\text{by the definition of }t_{x,y}\times\operatorname*{id}\nolimits_{M}%
\right) \nonumber\\
&  =\left(  t_{x,y}\left(  z\right)  ,q\right)  .
\label{pf.prop.sol.sigmacrosstau.transposX.cross1.LHS1}%
\end{align}

We have $q\in M=\left\{  1,2,\ldots,m\right\}  $. We are in one of the
following three cases:

\textit{Case 1:} We have $z=x$.

\textit{Case 2:} We have $z=y$.

\textit{Case 3:} We have neither $z=x$ nor $z=y$.

Let us first consider Case 1. In this case, we have $z=x$. But Lemma
\ref{lem.sol.perm.transX.tij1} \textbf{(a)} (applied to $x$ and $y$ instead of
$i$ and $j$) yields $t_{x,y}\left(  x\right)  =y$. But $t_{x,y}\left(
\underbrace{z}_{=x}\right)  =t_{x,y}\left(  x\right)  =y$.

Hence, (\ref{pf.prop.sol.sigmacrosstau.transposX.cross1.LHS1}) yields
\begin{equation}
\left(  t_{x,y}\times\operatorname*{id}\nolimits_{M}\right)  \left(  c\right)
=\left(  \underbrace{t_{x,y}\left(  z\right)  }_{=y},q\right)  =\left(
y,q\right)  . \label{pf.prop.sol.sigmacrosstau.transposX.cross1.c1.LHS}%
\end{equation}

On the other hand, $\left(  x,q\right)  \neq\left(  y,q\right)  $%
\ \ \ \ \footnote{\textit{Proof.} Assume the contrary. Thus, $\left(
x,q\right)  =\left(  y,q\right)  $. In other words, $x=y$ and $q=q$. But $x=y$
contradicts the fact that $x$ and $y$ are distinct. This contradiction shows
that our assumption was wrong. Qed.}. Hence, the two elements $\left(
x,q\right)  $ and $\left(  y,q\right)  $ of $X\times M$ are distinct. Thus,
Lemma \ref{lem.sol.perm.transX.tij1} \textbf{(a)} (applied to $X\times M$,
$\left(  x,q\right)  $ and $\left(  y,q\right)  $ instead of $X$, $i$ and $j$)
yields
\begin{equation}
t_{\left(  x,q\right)  ,\left(  y,q\right)  }\left(  \left(  x,q\right)
\right)  =\left(  y,q\right)  .
\label{pf.prop.sol.sigmacrosstau.transposX.cross1.c1.0}%
\end{equation}

Next, we observe that%
\begin{equation}
t_{\left(  x,j\right)  ,\left(  y,j\right)  }\left(  \left(  x,q\right)
\right)  =\left(  x,q\right)  \ \ \ \ \ \ \ \ \ \ \text{for each }j\in\left\{
1,2,\ldots,m\right\}  \text{ satisfying }j<q.
\label{pf.prop.sol.sigmacrosstau.transposX.cross1.c1.1}%
\end{equation}

[\textit{Proof of (\ref{pf.prop.sol.sigmacrosstau.transposX.cross1.c1.1}):}
Let $j\in\left\{  1,2,\ldots,m\right\}  $ be such that $j<q$. Thus, $j\neq q$
(since $j<q$). Hence, (\ref{pf.prop.sol.sigmacrosstau.transposX.cross1.t-dist}%
) (applied to $z=x$) yields $t_{\left(  x,j\right)  ,\left(  y,j\right)
}\left(  \left(  x,q\right)  \right)  =\left(  x,q\right)  $. This proves
(\ref{pf.prop.sol.sigmacrosstau.transposX.cross1.c1.1}).]

Next, we observe that%
\begin{equation}
t_{\left(  x,j\right)  ,\left(  y,j\right)  }\left(  \left(  y,q\right)
\right)  =\left(  y,q\right)  \ \ \ \ \ \ \ \ \ \ \text{for each }j\in\left\{
1,2,\ldots,m\right\}  \text{ satisfying }j>q.
\label{pf.prop.sol.sigmacrosstau.transposX.cross1.c1.2}%
\end{equation}

[\textit{Proof of (\ref{pf.prop.sol.sigmacrosstau.transposX.cross1.c1.2}):}
Let $j\in\left\{  1,2,\ldots,m\right\}  $ be such that $j>q$. Thus, $j\neq q$
(since $j>q$). Hence, (\ref{pf.prop.sol.sigmacrosstau.transposX.cross1.t-dist}%
) (applied to $z=y$) yields $t_{\left(  x,j\right)  ,\left(  y,j\right)
}\left(  \left(  y,q\right)  \right)  =\left(  y,q\right)  $. This proves
(\ref{pf.prop.sol.sigmacrosstau.transposX.cross1.c1.2}).]

Now, we can apply Lemma \ref{lem.sol.sigmacrosstau.compose1} to $X\times M$,
$t_{\left(  x,j\right)  ,\left(  y,j\right)  }$, $\left(  x,q\right)  $,
$\left(  y,q\right)  $ and $q$ instead of $X$, $f_{j}$, $x$, $y$ and $i$
(since the equalities (\ref{pf.prop.sol.sigmacrosstau.transposX.cross1.c1.0}),
(\ref{pf.prop.sol.sigmacrosstau.transposX.cross1.c1.1}) and
(\ref{pf.prop.sol.sigmacrosstau.transposX.cross1.c1.2}) hold). Thus, we obtain%
\begin{equation}
\left(  t_{\left(  x,m\right)  ,\left(  y,m\right)  }\circ t_{\left(
x,m-1\right)  ,\left(  y,m-1\right)  }\circ\cdots\circ t_{\left(  x,1\right)
,\left(  y,1\right)  }\right)  \left(  \left(  x,q\right)  \right)  =\left(
y,q\right)  . \label{pf.prop.sol.sigmacrosstau.transposX.cross1.c1.5}%
\end{equation}
In view of (\ref{pf.prop.sol.sigmacrosstau.transposX.cross1.s=}), this
rewrites as $s\left(  \left(  x,q\right)  \right)  =\left(  y,q\right)  $.
Comparing this with (\ref{pf.prop.sol.sigmacrosstau.transposX.cross1.c1.LHS}),
we obtain $\left(  t_{x,y}\times\operatorname*{id}\nolimits_{M}\right)
\left(  c\right)  =s\left(  \left(  x,q\right)  \right)  $. Comparing this
with $s\left(  \underbrace{c}_{=\left(  z,q\right)  }\right)  =s\left(
\left(  \underbrace{z}_{=x},q\right)  \right)  =s\left(  \left(  x,q\right)
\right)  $, we obtain $\left(  t_{x,y}\times\operatorname*{id}\nolimits_{M}%
\right)  \left(  c\right)  =s\left(  c\right)  $. Hence, $\left(
t_{x,y}\times\operatorname*{id}\nolimits_{M}\right)  \left(  c\right)
=s\left(  c\right)  $ is proven in Case 1.

Let us next consider Case 2. In this case, we have $z=y$. But Lemma
\ref{lem.sol.perm.transX.tij1} \textbf{(b)} (applied to $x$ and $y$ instead of
$i$ and $j$) yields $t_{x,y}\left(  y\right)  =x$. But $t_{x,y}\left(
\underbrace{z}_{=y}\right)  =t_{x,y}\left(  y\right)  =x$.

Hence, (\ref{pf.prop.sol.sigmacrosstau.transposX.cross1.LHS1}) yields
\begin{equation}
\left(  t_{x,y}\times\operatorname*{id}\nolimits_{M}\right)  \left(  c\right)
=\left(  \underbrace{t_{x,y}\left(  z\right)  }_{=x},q\right)  =\left(
x,q\right)  . \label{pf.prop.sol.sigmacrosstau.transposX.cross1.c2.LHS}%
\end{equation}

On the other hand, $\left(  x,q\right)  \neq\left(  y,q\right)  $%
\ \ \ \ \footnote{\textit{Proof.} Assume the contrary. Thus, $\left(
x,q\right)  =\left(  y,q\right)  $. In other words, $x=y$ and $q=q$. But $x=y$
contradicts the fact that $x$ and $y$ are distinct. This contradiction shows
that our assumption was wrong. Qed.}. Hence, the two elements $\left(
x,q\right)  $ and $\left(  y,q\right)  $ of $X\times M$ are distinct. Thus,
Lemma \ref{lem.sol.perm.transX.tij1} \textbf{(b)} (applied to $X\times M$,
$\left(  x,q\right)  $ and $\left(  y,q\right)  $ instead of $X$, $i$ and $j$)
yields
\begin{equation}
t_{\left(  x,q\right)  ,\left(  y,q\right)  }\left(  \left(  y,q\right)
\right)  =\left(  x,q\right)  .
\label{pf.prop.sol.sigmacrosstau.transposX.cross1.c2.0}%
\end{equation}

Next, we observe that%
\begin{equation}
t_{\left(  x,j\right)  ,\left(  y,j\right)  }\left(  \left(  y,q\right)
\right)  =\left(  y,q\right)  \ \ \ \ \ \ \ \ \ \ \text{for each }j\in\left\{
1,2,\ldots,m\right\}  \text{ satisfying }j<q.
\label{pf.prop.sol.sigmacrosstau.transposX.cross1.c2.1}%
\end{equation}

[\textit{Proof of (\ref{pf.prop.sol.sigmacrosstau.transposX.cross1.c2.1}):}
Let $j\in\left\{  1,2,\ldots,m\right\}  $ be such that $j<q$. Thus, $j\neq q$
(since $j<q$). Hence, (\ref{pf.prop.sol.sigmacrosstau.transposX.cross1.t-dist}%
) (applied to $z=y$) yields $t_{\left(  x,j\right)  ,\left(  y,j\right)
}\left(  \left(  y,q\right)  \right)  =\left(  y,q\right)  $. This proves
(\ref{pf.prop.sol.sigmacrosstau.transposX.cross1.c2.1}).]

Next, we observe that%
\begin{equation}
t_{\left(  x,j\right)  ,\left(  y,j\right)  }\left(  \left(  x,q\right)
\right)  =\left(  x,q\right)  \ \ \ \ \ \ \ \ \ \ \text{for each }j\in\left\{
1,2,\ldots,m\right\}  \text{ satisfying }j>q.
\label{pf.prop.sol.sigmacrosstau.transposX.cross1.c2.2}%
\end{equation}

[\textit{Proof of (\ref{pf.prop.sol.sigmacrosstau.transposX.cross1.c2.2}):}
Let $j\in\left\{  1,2,\ldots,m\right\}  $ be such that $j>q$. Thus, $j\neq q$
(since $j>q$). Hence, (\ref{pf.prop.sol.sigmacrosstau.transposX.cross1.t-dist}%
) (applied to $z=x$) yields $t_{\left(  x,j\right)  ,\left(  y,j\right)
}\left(  \left(  x,q\right)  \right)  =\left(  x,q\right)  $. This proves
(\ref{pf.prop.sol.sigmacrosstau.transposX.cross1.c2.2}).]

Now, we can apply Lemma \ref{lem.sol.sigmacrosstau.compose1} to $X\times M$,
$t_{\left(  x,j\right)  ,\left(  y,j\right)  }$, $\left(  y,q\right)  $,
$\left(  x,q\right)  $ and $q$ instead of $X$, $f_{j}$, $x$, $y$ and $i$
(since the equalities (\ref{pf.prop.sol.sigmacrosstau.transposX.cross1.c2.0}),
(\ref{pf.prop.sol.sigmacrosstau.transposX.cross1.c2.1}) and
(\ref{pf.prop.sol.sigmacrosstau.transposX.cross1.c2.2}) hold). Thus, we obtain%
\begin{equation}
\left(  t_{\left(  x,m\right)  ,\left(  y,m\right)  }\circ t_{\left(
x,m-1\right)  ,\left(  y,m-1\right)  }\circ\cdots\circ t_{\left(  x,1\right)
,\left(  y,1\right)  }\right)  \left(  \left(  y,q\right)  \right)  =\left(
x,q\right)  . \label{pf.prop.sol.sigmacrosstau.transposX.cross1.c2.5}%
\end{equation}
In view of (\ref{pf.prop.sol.sigmacrosstau.transposX.cross1.s=}), this
rewrites as $s\left(  \left(  y,q\right)  \right)  =\left(  x,q\right)  $.
Comparing this with (\ref{pf.prop.sol.sigmacrosstau.transposX.cross1.c2.LHS}),
we obtain $\left(  t_{x,y}\times\operatorname*{id}\nolimits_{M}\right)
\left(  c\right)  =s\left(  \left(  y,q\right)  \right)  $. Comparing this
with $s\left(  \underbrace{c}_{=\left(  z,q\right)  }\right)  =s\left(
\left(  \underbrace{z}_{=y},q\right)  \right)  =s\left(  \left(  y,q\right)
\right)  $, we obtain $\left(  t_{x,y}\times\operatorname*{id}\nolimits_{M}%
\right)  \left(  c\right)  =s\left(  c\right)  $. Hence, $\left(
t_{x,y}\times\operatorname*{id}\nolimits_{M}\right)  \left(  c\right)
=s\left(  c\right)  $ is proven in Case 2.

Let us finally consider Case 3. In this case, we have neither $z=x$ nor $z=y$.
Thus, $z\notin\left\{  x,y\right\}  $\ \ \ \ \footnote{\textit{Proof.} Assume
the contrary. Thus, $z\in\left\{  x,y\right\}  $. In other words, $z$ equals
either $x$ or $y$. In other words, we have either $z=x$ or $z=y$. But this
contradicts the fact that we have neither $z=x$ nor $z=y$. This contradiction
shows that our assumption was wrong. Qed.}. Combining this with $z\in X$, we
obtain $z\in X\setminus\left\{  x,y\right\}  $. Thus, Lemma
\ref{lem.sol.perm.transX.tij1} \textbf{(c)} (applied to $x$, $y$ and $z$
instead of $i$, $j$ and $k$) yields $t_{x,y}\left(  z\right)  =z$.

Hence, (\ref{pf.prop.sol.sigmacrosstau.transposX.cross1.LHS1}) yields
\begin{equation}
\left(  t_{x,y}\times\operatorname*{id}\nolimits_{M}\right)  \left(  c\right)
=\left(  \underbrace{t_{x,y}\left(  z\right)  }_{=z},q\right)  =\left(
z,q\right)  . \label{pf.prop.sol.sigmacrosstau.transposX.cross1.c3.LHS}%
\end{equation}

Next, we observe that%
\begin{equation}
t_{\left(  x,j\right)  ,\left(  y,j\right)  }\left(  \left(  z,q\right)
\right)  =\left(  z,q\right)  \ \ \ \ \ \ \ \ \ \ \text{for each }j\in\left\{
1,2,\ldots,m\right\}  .
\label{pf.prop.sol.sigmacrosstau.transposX.cross1.c3.1}%
\end{equation}

[\textit{Proof of (\ref{pf.prop.sol.sigmacrosstau.transposX.cross1.c3.1}):}
Let $j\in\left\{  1,2,\ldots,m\right\}  $. Thus, $j\in\left\{  1,2,\ldots
,m\right\}  =M$. The elements $\left(  x,j\right)  $ and $\left(  y,j\right)
$ of $X\times M$ are distinct\footnote{\textit{Proof.} Assume the contrary.
Thus, $\left(  x,j\right)  =\left(  y,j\right)  $. In other words, $x=y$ and
$j=j$. But $x=y$ contradicts the fact that $x$ and $y$ are distinct. This
contradiction shows that our assumption was wrong. Qed.}.

We have $\left(  z,q\right)  \neq\left(  x,j\right)  $%
\ \ \ \ \footnote{\textit{Proof.} Assume the contrary. Thus, $\left(
z,q\right)  =\left(  x,j\right)  $. In other words, $z=x$ and $q=j$. But $z=x$
contradicts the fact that we have neither $z=x$ nor $z=y$. This contradiction
shows that our assumption was false. Qed.}. In other words, we don't have
$\left(  z,q\right)  =\left(  x,j\right)  $.

Also, $\left(  z,q\right)  \neq\left(  y,j\right)  $%
\ \ \ \ \footnote{\textit{Proof.} Assume the contrary. Thus, $\left(
z,q\right)  =\left(  y,j\right)  $. In other words, $z=y$ and $q=j$. But $z=y$
contradicts the fact that we have neither $z=x$ nor $z=y$. This contradiction
shows that our assumption was false. Qed.}. In other words, we don't have
$\left(  z,q\right)  =\left(  y,j\right)  $. Hence, we have neither $\left(
z,q\right)  =\left(  x,j\right)  $ nor $\left(  z,q\right)  =\left(
y,j\right)  $.

We have $\left(  z,q\right)  \notin\left\{  \left(  x,j\right)  ,\left(
y,j\right)  \right\}  $\ \ \ \ \footnote{\textit{Proof.} Assume the contrary.
Thus, $\left(  z,q\right)  \in\left\{  \left(  x,j\right)  ,\left(
y,j\right)  \right\}  $. In other words, $\left(  z,q\right)  $ equals either
$\left(  x,j\right)  $ or $\left(  y,j\right)  $. In other words, we have
either $\left(  z,q\right)  =\left(  x,j\right)  $ or $\left(  z,q\right)
=\left(  y,j\right)  $. This contradicts the fact that we have neither
$\left(  z,q\right)  =\left(  x,j\right)  $ nor $\left(  z,q\right)  =\left(
y,j\right)  $. This contradiction shows that our assumption was wrong. Qed.}.
Combining $\left(  z,q\right)  \in X\times M$ (since $z\in X$ and $q\in M$)
with $\left(  z,q\right)  \notin\left\{  \left(  x,j\right)  ,\left(
y,j\right)  \right\}  $, we obtain $\left(  z,q\right)  \in\left(  X\times
M\right)  \setminus\left\{  \left(  x,j\right)  ,\left(  y,j\right)  \right\}
$. Hence, Lemma \ref{lem.sol.perm.transX.tij1} \textbf{(c)} (applied to
$X\times M$, $\left(  x,j\right)  $, $\left(  y,j\right)  $ and $\left(
z,q\right)  $ instead of $X$, $i$, $j$ and $k$) yields $t_{\left(  x,j\right)
,\left(  y,j\right)  }\left(  \left(  z,q\right)  \right)  =\left(
z,q\right)  $. This proves
(\ref{pf.prop.sol.sigmacrosstau.transposX.cross1.c3.1}).]

Now, we can apply Lemma \ref{lem.sol.sigmacrosstau.compose0} to $X\times M$,
$t_{\left(  x,j\right)  ,\left(  y,j\right)  }$ and $\left(  z,q\right)  $
instead of $X$, $f_{j}$ and $x$ (since
(\ref{pf.prop.sol.sigmacrosstau.transposX.cross1.c3.1}) holds). Thus, we
obtain%
\begin{equation}
\left(  t_{\left(  x,m\right)  ,\left(  y,m\right)  }\circ t_{\left(
x,m-1\right)  ,\left(  y,m-1\right)  }\circ\cdots\circ t_{\left(  x,1\right)
,\left(  y,1\right)  }\right)  \left(  \left(  z,q\right)  \right)  =\left(
z,q\right)  . \label{pf.prop.sol.sigmacrosstau.transposX.cross1.c3.5}%
\end{equation}
In view of (\ref{pf.prop.sol.sigmacrosstau.transposX.cross1.s=}), this
rewrites as $s\left(  \left(  z,q\right)  \right)  =\left(  z,q\right)  $.
Comparing this with (\ref{pf.prop.sol.sigmacrosstau.transposX.cross1.c3.LHS}),
we obtain $\left(  t_{x,y}\times\operatorname*{id}\nolimits_{M}\right)
\left(  c\right)  =s\left(  \left(  z,q\right)  \right)  $. Comparing this
with $s\left(  \underbrace{c}_{=\left(  z,q\right)  }\right)  =s\left(
\left(  z,q\right)  \right)  $, we obtain $\left(  t_{x,y}\times
\operatorname*{id}\nolimits_{M}\right)  \left(  c\right)  =s\left(  c\right)
$. Hence, $\left(  t_{x,y}\times\operatorname*{id}\nolimits_{M}\right)
\left(  c\right)  =s\left(  c\right)  $ is proven in Case 3.

We have now proven $\left(  t_{x,y}\times\operatorname*{id}\nolimits_{M}%
\right)  \left(  c\right)  =s\left(  c\right)  $ in each of the three Cases 1,
2 and 3. Since these three Cases cover all possibilities, this shows that
$\left(  t_{x,y}\times\operatorname*{id}\nolimits_{M}\right)  \left(
c\right)  =s\left(  c\right)  $ always holds.

Now, forget that we fixed $c$. We thus have proven that $\left(  t_{x,y}%
\times\operatorname*{id}\nolimits_{M}\right)  \left(  c\right)  =s\left(
c\right)  $ for each $c\in X\times M$. Thus,%
\[
t_{x,y}\times\operatorname*{id}\nolimits_{M}=s=t_{\left(  x,m\right)  ,\left(
y,m\right)  }\circ t_{\left(  x,m-1\right)  ,\left(  y,m-1\right)  }%
\circ\cdots\circ t_{\left(  x,1\right)  ,\left(  y,1\right)  }.
\]
This completes the proof of Proposition
\ref{prop.sol.sigmacrosstau.transposX.cross1}.
\end{proof}
\end{verlong}

\begin{corollary}
\label{cor.sol.sigmacrosstau.transposX.cross2}Let $X$ be a finite set. Let
$m\in\mathbb{N}$. Let $M$ be the set $\left\{  1,2,\ldots,m\right\}  $. Let
$x$ and $y$ be two distinct elements of $X$. Then, the permutation
$t_{x,y}\times\operatorname*{id}\nolimits_{M}$ of $X\times M$ satisfies
$\left(  -1\right)  ^{t_{x,y}\times\operatorname*{id}\nolimits_{M}}=\left(
-1\right)  ^{m}$.
\end{corollary}

\begin{vershort}
\begin{proof}
[Proof of Corollary \ref{cor.sol.sigmacrosstau.transposX.cross2}.]Corollary
\ref{cor.sol.sigmacrosstau.permxperm} (applied to $X$, $M$, $t_{x,y}$ and
$\operatorname*{id}\nolimits_{M}$ instead of $U$, $V$, $\sigma$ and $\tau$)
shows that $t_{x,y}\times\operatorname*{id}\nolimits_{M}$ is a permutation of
$X\times M$. Proposition \ref{prop.sol.sigmacrosstau.transposX.cross1} shows
that $t_{\left(  x,1\right)  ,\left(  y,1\right)  },t_{\left(  x,2\right)
,\left(  y,2\right)  },\ldots,t_{\left(  x,m\right)  ,\left(  y,m\right)  }$
are $m$ well-defined transpositions of the set $X\times M$, and that we have%
\begin{equation}
t_{x,y}\times\operatorname*{id}\nolimits_{M}=t_{\left(  x,m\right)  ,\left(
y,m\right)  }\circ t_{\left(  x,m-1\right)  ,\left(  y,m-1\right)  }%
\circ\cdots\circ t_{\left(  x,1\right)  ,\left(  y,1\right)  }.
\label{pf.cor.sol.sigmacrosstau.transposX.cross2.short.2}%
\end{equation}
The equality (\ref{pf.cor.sol.sigmacrosstau.transposX.cross2.short.2}) shows
that $t_{x,y}\times\operatorname*{id}\nolimits_{M}$ can be written as a
composition of $m$ transpositions of $X\times M$ (since $t_{\left(
x,1\right)  ,\left(  y,1\right)  },t_{\left(  x,2\right)  ,\left(  y,2\right)
},\ldots,t_{\left(  x,m\right)  ,\left(  y,m\right)  }$ are $m$ transpositions
of $X\times M$). Hence, Corollary \ref{cor.sol.perm.transX.sign.prod-of-trans}
(applied to $X\times M$, $t_{x,y}\times\operatorname*{id}\nolimits_{M}$ and
$m$ instead of $X$, $\sigma$ and $k$) shows that $\left(  -1\right)
^{t_{x,y}\times\operatorname*{id}\nolimits_{M}}=\left(  -1\right)  ^{m}$. This
proves Corollary \ref{cor.sol.sigmacrosstau.transposX.cross2}.
\end{proof}
\end{vershort}

\begin{verlong}
\begin{proof}
[Proof of Corollary \ref{cor.sol.sigmacrosstau.transposX.cross2}.]The sets $X$
and $M$ are finite. Hence, their product $X\times M$ is also finite.

Also, $t_{x,y}$ is a permutation of $X$, whereas $\operatorname*{id}%
\nolimits_{M}$ is a permutation of $M$. Hence, $t_{x,y}\times
\operatorname*{id}\nolimits_{M}$ is a permutation of $X\times M$ (by Corollary
\ref{cor.sol.sigmacrosstau.permxperm} (applied to $X$, $M$, $t_{x,y}$ and
$\operatorname*{id}\nolimits_{M}$ instead of $U$, $V$, $\sigma$ and $\tau$)).

Proposition \ref{prop.sol.sigmacrosstau.transposX.cross1} shows that
$t_{\left(  x,1\right)  ,\left(  y,1\right)  },t_{\left(  x,2\right)  ,\left(
y,2\right)  },\ldots,t_{\left(  x,m\right)  ,\left(  y,m\right)  }$ are $m$
well-defined transpositions of the set $X\times M$, and that we have%
\begin{equation}
t_{x,y}\times\operatorname*{id}\nolimits_{M}=t_{\left(  x,m\right)  ,\left(
y,m\right)  }\circ t_{\left(  x,m-1\right)  ,\left(  y,m-1\right)  }%
\circ\cdots\circ t_{\left(  x,1\right)  ,\left(  y,1\right)  }.
\label{pf.cor.sol.sigmacrosstau.transposX.cross2.2}%
\end{equation}

Now, we know that $t_{\left(  x,1\right)  ,\left(  y,1\right)  },t_{\left(
x,2\right)  ,\left(  y,2\right)  },\ldots,t_{\left(  x,m\right)  ,\left(
y,m\right)  }$ are $m$ transpositions of $X\times M$. In other words,
$t_{\left(  x,m\right)  ,\left(  y,m\right)  },t_{\left(  x,m-1\right)
,\left(  y,m-1\right)  },\ldots,t_{\left(  x,1\right)  ,\left(  y,1\right)  }$
are $m$ transpositions of $X\times M$. Hence, $t_{\left(  x,m\right)  ,\left(
y,m\right)  }\circ t_{\left(  x,m-1\right)  ,\left(  y,m-1\right)  }%
\circ\cdots\circ t_{\left(  x,1\right)  ,\left(  y,1\right)  }$ is a
composition of $m$ transpositions of $X\times M$. Thus, the equality
(\ref{pf.cor.sol.sigmacrosstau.transposX.cross2.2}) shows that $t_{x,y}%
\times\operatorname*{id}\nolimits_{M}$ can be written as a composition of $m$
transpositions of $X\times M$. Hence, Corollary
\ref{cor.sol.perm.transX.sign.prod-of-trans} (applied to $X\times M$,
$t_{x,y}\times\operatorname*{id}\nolimits_{M}$ and $m$ instead of $X$,
$\sigma$ and $k$) shows that $\left(  -1\right)  ^{t_{x,y}\times
\operatorname*{id}\nolimits_{M}}=\left(  -1\right)  ^{m}$. This proves
Corollary \ref{cor.sol.sigmacrosstau.transposX.cross2}.
\end{proof}
\end{verlong}

\begin{corollary}
\label{cor.sol.sigmacrosstau.transposX.cross3}Let $U$ and $V$ be two finite
sets. Let $i$ and $j$ be two distinct elements of $U$. Then, the permutation
$t_{i,j}\times\operatorname*{id}\nolimits_{V}$ of $U\times V$ satisfies
$\left(  -1\right)  ^{t_{i,j}\times\operatorname*{id}\nolimits_{V}}=\left(
-1\right)  ^{\left\vert V\right\vert }$.
\end{corollary}

\begin{proof}
[Proof of Corollary \ref{cor.sol.sigmacrosstau.transposX.cross3}.]Define an
$m\in\mathbb{N}$ by $m=\left\vert V\right\vert $. (This is well-defined, since
the set $V$ is finite.)

Let $M$ be the set $\left\{  1,2,\ldots,m\right\}  $. Thus, $M=\left\{
1,2,\ldots,m\right\}  $, so that $\left\vert M\right\vert =\left\vert \left\{
1,2,\ldots,m\right\}  \right\vert =m=\left\vert V\right\vert $. The sets $V$
and $M$ are finite and have the same size (since $\left\vert M\right\vert
=\left\vert V\right\vert $). Hence, there exists a bijection $\phi
:M\rightarrow V$. Consider such a $\phi$.

\begin{vershort}
Corollary \ref{cor.sol.sigmacrosstau.bijxbij} (applied to $U$, $U$, $M$, $V$,
$\operatorname*{id}\nolimits_{U}$ and $\phi$ instead of $X$, $X^{\prime}$,
$Y$, $Y^{\prime}$, $\alpha$ and $\beta$) shows that the map
$\operatorname*{id}\nolimits_{U}\times\phi:U\times M\rightarrow U\times V$ is
bijective. In other words, $\operatorname*{id}\nolimits_{U}\times\phi$ is a
bijection. Also, Corollary \ref{cor.sol.sigmacrosstau.permxperm} (applied to
$M$, $t_{i,j}$ and $\operatorname*{id}\nolimits_{M}$ instead of $V$, $\sigma$
and $\tau$) shows that $t_{i,j}\times\operatorname*{id}\nolimits_{M}$ is a
permutation of $U\times M$. Similarly, $t_{i,j}\times\operatorname*{id}%
\nolimits_{V}$ is a permutation of $U\times V$.
\end{vershort}

\begin{verlong}
The set $U\times M$ is finite (since it is the product of the two finite sets
$U$ and $M$). The set $U\times V$ is finite (since it is the product of the
two finite sets $U$ and $V$).

The maps $\operatorname*{id}\nolimits_{U}:U\rightarrow U$ and $\phi
:M\rightarrow V$ are bijections, i.e., bijective maps. Thus, Corollary
\ref{cor.sol.sigmacrosstau.bijxbij} (applied to $U$, $U$, $M$, $V$,
$\operatorname*{id}\nolimits_{U}$ and $\phi$ instead of $X$, $X^{\prime}$,
$Y$, $Y^{\prime}$, $\alpha$ and $\beta$) shows that the map
$\operatorname*{id}\nolimits_{U}\times\phi:U\times M\rightarrow U\times V$ is
bijective as well, and its inverse is the map $\left(  \operatorname*{id}%
\nolimits_{U}\right)  ^{-1}\times\phi^{-1}:U\times V\rightarrow U\times M$. In
particular, $\operatorname*{id}\nolimits_{U}\times\phi:U\times M\rightarrow
U\times V$ is bijective. In other words, $\operatorname*{id}\nolimits_{U}%
\times\phi:U\times M\rightarrow U\times V$ is a bijection.

Also, $t_{i,j}$ is a permutation of $U$ (since $i$ and $j$ are two distinct
elements of $U$), whereas $\operatorname*{id}\nolimits_{M}$ is a permutation
of $M$. Hence, $t_{i,j}\times\operatorname*{id}\nolimits_{M}$ is a permutation
of $U\times M$ (by Corollary \ref{cor.sol.sigmacrosstau.permxperm} (applied to
$M$, $t_{i,j}$ and $\operatorname*{id}\nolimits_{M}$ instead of $V$, $\sigma$
and $\tau$)).

Also, $t_{i,j}$ is a permutation of $U$, whereas $\operatorname*{id}%
\nolimits_{V}$ is a permutation of $V$. Hence, $t_{i,j}\times
\operatorname*{id}\nolimits_{V}$ is a permutation of $U\times V$ (by Corollary
\ref{cor.sol.sigmacrosstau.permxperm} (applied to $t_{i,j}$ and
$\operatorname*{id}\nolimits_{V}$ instead of $\sigma$ and $\tau$)).
\end{verlong}

Now, Proposition \ref{prop.sol.sigmacrosstau.signs-conj} (applied to $U\times
M$, $U\times V$, $\operatorname*{id}\nolimits_{U}\times\phi$ and
$t_{i,j}\times\operatorname*{id}\nolimits_{M}$ instead of $X$, $Y$, $f$ and
$\sigma$) shows that $\left(  \operatorname*{id}\nolimits_{U}\times
\phi\right)  \circ\left(  t_{i,j}\times\operatorname*{id}\nolimits_{M}\right)
\circ\left(  \operatorname*{id}\nolimits_{U}\times\phi\right)  ^{-1}$ is a
permutation of $U\times V$ and satisfies
\begin{equation}
\left(  -1\right)  ^{\left(  \operatorname*{id}\nolimits_{U}\times\phi\right)
\circ\left(  t_{i,j}\times\operatorname*{id}\nolimits_{M}\right)  \circ\left(
\operatorname*{id}\nolimits_{U}\times\phi\right)  ^{-1}}=\left(  -1\right)
^{t_{i,j}\times\operatorname*{id}\nolimits_{M}}.
\label{pf.cor.sol.sigmacrosstau.transposX.cross3.3}%
\end{equation}
But Corollary \ref{cor.sol.sigmacrosstau.transposX.cross2} (applied to $X=U$,
$x=i$ and $y=j$) shows that the permutation $t_{i,j}\times\operatorname*{id}%
\nolimits_{M}$ of $U\times M$ satisfies%
\begin{equation}
\left(  -1\right)  ^{t_{i,j}\times\operatorname*{id}\nolimits_{M}}=\left(
-1\right)  ^{m}. \label{pf.cor.sol.sigmacrosstau.transposX.cross3.4}%
\end{equation}

But Proposition \ref{prop.sol.sigmacrosstau.fgf'g'} (applied to $U$, $U$, $U$,
$M$, $M$, $V$, $t_{i,j}$, $\operatorname*{id}\nolimits_{U}$,
$\operatorname*{id}\nolimits_{M}$ and $\phi$ instead of $X$, $X^{\prime}$,
$X^{\prime\prime}$, $Y$, $Y^{\prime}$, $Y^{\prime\prime}$, $\alpha$,
$\alpha^{\prime}$, $\beta$ and $\beta^{\prime}$) yields%
\begin{equation}
\left(  \operatorname*{id}\nolimits_{U}\times\phi\right)  \circ\left(
t_{i,j}\times\operatorname*{id}\nolimits_{M}\right)  =\underbrace{\left(
\operatorname*{id}\nolimits_{U}\circ t_{i,j}\right)  }_{=t_{i,j}}%
\times\underbrace{\left(  \phi\circ\operatorname*{id}\nolimits_{M}\right)
}_{=\phi}=t_{i,j}\times\phi.
\label{pf.cor.sol.sigmacrosstau.transposX.cross3.5}%
\end{equation}
On the other hand, Proposition \ref{prop.sol.sigmacrosstau.fgf'g'} (applied to
$U$, $U$, $U$, $M$, $V$, $V$, $\operatorname*{id}\nolimits_{U}$, $t_{i,j}$,
$\phi$ and $\operatorname*{id}\nolimits_{V}$ instead of $X$, $X^{\prime}$,
$X^{\prime\prime}$, $Y$, $Y^{\prime}$, $Y^{\prime\prime}$, $\alpha$,
$\alpha^{\prime}$, $\beta$ and $\beta^{\prime}$) yields%
\[
\left(  t_{i,j}\times\operatorname*{id}\nolimits_{V}\right)  \circ\left(
\operatorname*{id}\nolimits_{U}\times\phi\right)  =\underbrace{\left(
t_{i,j}\circ\operatorname*{id}\nolimits_{U}\right)  }_{=t_{i,j}}%
\times\underbrace{\left(  \operatorname*{id}\nolimits_{V}\circ\phi\right)
}_{=\phi}=t_{i,j}\times\phi.
\]
Comparing this with (\ref{pf.cor.sol.sigmacrosstau.transposX.cross3.5}), we
obtain%
\[
\left(  \operatorname*{id}\nolimits_{U}\times\phi\right)  \circ\left(
t_{i,j}\times\operatorname*{id}\nolimits_{M}\right)  =\left(  t_{i,j}%
\times\operatorname*{id}\nolimits_{V}\right)  \circ\left(  \operatorname*{id}%
\nolimits_{U}\times\phi\right)  .
\]
Hence,%
\begin{align*}
&  \underbrace{\left(  \operatorname*{id}\nolimits_{U}\times\phi\right)
\circ\left(  t_{i,j}\times\operatorname*{id}\nolimits_{M}\right)  }_{=\left(
t_{i,j}\times\operatorname*{id}\nolimits_{V}\right)  \circ\left(
\operatorname*{id}\nolimits_{U}\times\phi\right)  }\circ\left(
\operatorname*{id}\nolimits_{U}\times\phi\right)  ^{-1}\\
&  =\left(  t_{i,j}\times\operatorname*{id}\nolimits_{V}\right)
\circ\underbrace{\left(  \operatorname*{id}\nolimits_{U}\times\phi\right)
\circ\left(  \operatorname*{id}\nolimits_{U}\times\phi\right)  ^{-1}%
}_{=\operatorname*{id}\nolimits_{U\times V}}=t_{i,j}\times\operatorname*{id}%
\nolimits_{V}.
\end{align*}
Thus, (\ref{pf.cor.sol.sigmacrosstau.transposX.cross3.3}) rewrites as $\left(
-1\right)  ^{t_{i,j}\times\operatorname*{id}\nolimits_{V}}=\left(  -1\right)
^{t_{i,j}\times\operatorname*{id}\nolimits_{M}}$. Thus,%
\begin{align*}
\left(  -1\right)  ^{t_{i,j}\times\operatorname*{id}\nolimits_{V}}  &
=\left(  -1\right)  ^{t_{i,j}\times\operatorname*{id}\nolimits_{M}}=\left(
-1\right)  ^{m}\ \ \ \ \ \ \ \ \ \ \left(  \text{by
(\ref{pf.cor.sol.sigmacrosstau.transposX.cross3.4})}\right) \\
&  =\left(  -1\right)  ^{\left\vert V\right\vert }\ \ \ \ \ \ \ \ \ \ \left(
\text{since }m=\left\vert V\right\vert \right)  .
\end{align*}
This proves Corollary \ref{cor.sol.sigmacrosstau.transposX.cross3}.
\end{proof}

\begin{corollary}
\label{cor.sol.sigmacrosstau.sign-sigma-cross-id}Let $U$ and $V$ be two finite
sets. Let $\sigma$ be a permutation of $U$. Then, the permutation
$\sigma\times\operatorname*{id}\nolimits_{V}$ of $U\times V$ satisfies
$\left(  -1\right)  ^{\sigma\times\operatorname*{id}\nolimits_{V}}=\left(
\left(  -1\right)  ^{\sigma}\right)  ^{\left\vert V\right\vert }$.
\end{corollary}

\begin{vershort}
\begin{proof}
[Proof of Corollary \ref{cor.sol.sigmacrosstau.sign-sigma-cross-id}.]Corollary
\ref{cor.sol.sigmacrosstau.permxperm} (applied to $\tau=\operatorname*{id}%
\nolimits_{V}$) shows that $\sigma\times\operatorname*{id}\nolimits_{V}$ is a
permutation of $U\times V$.

Proposition \ref{prop.sol.perm.transX.product} (applied to $U$ and $\sigma$
instead of $X$ and $\tau$) shows that $\sigma$ can be written as a composition
of finitely many transpositions of $U$. In other words, there exists some
$k\in\mathbb{N}$ such that $\sigma$ can be written as a composition of $k$
transpositions of $U$. Consider this $k$. Corollary
\ref{cor.sol.perm.transX.sign.prod-of-trans} (applied to $X=U$) yields%
\begin{equation}
\left(  -1\right)  ^{\sigma}=\left(  -1\right)  ^{k}.
\label{pf.cor.sol.sigmacrosstau.sign-sigma-cross-id.short.1}%
\end{equation}

We know that $\sigma$ can be written as a composition of $k$ transpositions of
$U$. In other words, there exist $k$ transpositions $f_{1},f_{2},\ldots,f_{k}$
of $U$ satisfying $\sigma=f_{1}\circ f_{2}\circ\cdots\circ f_{k}$. Consider
these $f_{1},f_{2},\ldots,f_{k}$.

For each $i\in\left\{  1,2,\ldots,k\right\}  $, define a map $g_{i}:U\times
V\rightarrow U\times V$ by $g_{i}=f_{i}\times\operatorname*{id}\nolimits_{V}$.

Let $p\in\left\{  1,2,\ldots,k\right\}  $. Thus, $g_{p}=f_{p}\times
\operatorname*{id}\nolimits_{V}$ (by the definition of $g_{p}$). Corollary
\ref{cor.sol.sigmacrosstau.permxperm} (applied to $f_{p}$ and
$\operatorname*{id}\nolimits_{V}$ instead of $\sigma$ and $\tau$) shows that
$f_{p}\times\operatorname*{id}\nolimits_{V}$ is a permutation of $U\times V$.
In other words, $g_{p}$ is a permutation of $U\times V$ (since $g_{p}%
=f_{p}\times\operatorname*{id}\nolimits_{V}$).

The map $f_{p}$ is a transposition of $U$. In other words, $f_{p}=t_{i,j}$ for
some two distinct elements $i$ and $j$ of $U$ (by the definition of a
transposition of $U$). Consider these $i$ and $j$. Corollary
\ref{cor.sol.sigmacrosstau.transposX.cross3} shows that the permutation
$t_{i,j}\times\operatorname*{id}\nolimits_{V}$ of $U\times V$ satisfies
$\left(  -1\right)  ^{t_{i,j}\times\operatorname*{id}\nolimits_{V}}=\left(
-1\right)  ^{\left\vert V\right\vert }$. Now, $f_{p}=t_{i,j}$. Thus,
$f_{p}\times\operatorname*{id}\nolimits_{V}=t_{i,j}\times\operatorname*{id}%
\nolimits_{V}$. Hence, $g_{p}=f_{p}\times\operatorname*{id}\nolimits_{V}%
=t_{i,j}\times\operatorname*{id}\nolimits_{V}$. Therefore,%
\begin{equation}
\left(  -1\right)  ^{g_{p}}=\left(  -1\right)  ^{t_{i,j}\times
\operatorname*{id}\nolimits_{V}}=\left(  -1\right)  ^{\left\vert V\right\vert
}. \label{pf.cor.sol.sigmacrosstau.sign-sigma-cross-id.short.3}%
\end{equation}

Now, forget that we fixed $p$. We thus have proven that for each $p\in\left\{
1,2,\ldots,k\right\}  $, the map $g_{p}$ is a permutation of $U\times V$ and
satisfies (\ref{pf.cor.sol.sigmacrosstau.sign-sigma-cross-id.short.3}).

Hence, the maps $g_{1},g_{2},\ldots,g_{k}$ are $k$ permutations of $U\times
V$. Thus, Proposition \ref{prop.sol.perm.transX.sign.prod-of-many} (applied to
$X=U\times V$ and $\sigma_{h}=g_{h}$) yields%
\begin{equation}
\left(  -1\right)  ^{g_{1}\circ g_{2}\circ\cdots\circ g_{k}}=\left(
-1\right)  ^{g_{1}}\cdot\left(  -1\right)  ^{g_{2}}\cdot\cdots\cdot\left(
-1\right)  ^{g_{k}}.
\label{pf.cor.sol.sigmacrosstau.sign-sigma-cross-id.short.5}%
\end{equation}

Corollary \ref{cor.sol.sigmacrosstau.fixid} yields%
\[
g_{1}\circ g_{2}\circ\cdots\circ g_{k}=\underbrace{\left(  f_{1}\circ
f_{2}\circ\cdots\circ f_{k}\right)  }_{\substack{=\sigma\\\text{(since }%
\sigma=f_{1}\circ f_{2}\circ\cdots\circ f_{k}\text{)}}}\times
\operatorname*{id}\nolimits_{V}=\sigma\times\operatorname*{id}\nolimits_{V}.
\]
Thus, $\sigma\times\operatorname*{id}\nolimits_{V}=g_{1}\circ g_{2}\circ
\cdots\circ g_{k}$. Therefore,%
\begin{align*}
\left(  -1\right)  ^{\sigma\times\operatorname*{id}\nolimits_{V}}  &  =\left(
-1\right)  ^{g_{1}\circ g_{2}\circ\cdots\circ g_{k}}=\left(  -1\right)
^{g_{1}}\cdot\left(  -1\right)  ^{g_{2}}\cdot\cdots\cdot\left(  -1\right)
^{g_{k}}\\
&  =\prod_{p=1}^{k}\underbrace{\left(  -1\right)  ^{g_{p}}}%
_{\substack{=\left(  -1\right)  ^{\left\vert V\right\vert }\\\text{(by
(\ref{pf.cor.sol.sigmacrosstau.sign-sigma-cross-id.short.3}))}}}=\prod
_{p=1}^{k}\left(  -1\right)  ^{\left\vert V\right\vert }=\left(  \left(
-1\right)  ^{\left\vert V\right\vert }\right)  ^{k}\\
&  =\left(  -1\right)  ^{\left\vert V\right\vert \cdot k}=\left(  -1\right)
^{k\cdot\left\vert V\right\vert }=\left(  \underbrace{\left(  -1\right)  ^{k}%
}_{\substack{=\left(  -1\right)  ^{\sigma}\\\text{(by
(\ref{pf.cor.sol.sigmacrosstau.sign-sigma-cross-id.short.1}))}}}\right)
^{\left\vert V\right\vert }=\left(  \left(  -1\right)  ^{\sigma}\right)
^{\left\vert V\right\vert }.
\end{align*}
This proves Corollary \ref{cor.sol.sigmacrosstau.sign-sigma-cross-id}.
\end{proof}
\end{vershort}

\begin{verlong}
\begin{proof}
[Proof of Corollary \ref{cor.sol.sigmacrosstau.sign-sigma-cross-id}.]The map
$\sigma$ is a permutation of $U$, whereas $\operatorname*{id}\nolimits_{V}$ is
a permutation of $V$. Hence, $\sigma\times\operatorname*{id}\nolimits_{V}$ is
a permutation of $U\times V$ (by Corollary
\ref{cor.sol.sigmacrosstau.permxperm} (applied to $\tau=\operatorname*{id}%
\nolimits_{V}$)).

Proposition \ref{prop.sol.perm.transX.product} (applied to $U$ and $\sigma$
instead of $X$ and $\tau$) shows that $\sigma$ can be written as a composition
of finitely many transpositions of $U$. In other words, there exists some
$k\in\mathbb{N}$ such that $\sigma$ can be written as a composition of $k$
transpositions of $U$. Consider this $k$. Corollary
\ref{cor.sol.perm.transX.sign.prod-of-trans} (applied to $X=U$) yields%
\begin{equation}
\left(  -1\right)  ^{\sigma}=\left(  -1\right)  ^{k}.
\label{pf.cor.sol.sigmacrosstau.sign-sigma-cross-id.1}%
\end{equation}

We know that $\sigma$ can be written as a composition of $k$ transpositions of
$U$. In other words, there exist $k$ transpositions $f_{1},f_{2},\ldots,f_{k}$
of $U$ satisfying $\sigma=f_{1}\circ f_{2}\circ\cdots\circ f_{k}$. Consider
these $f_{1},f_{2},\ldots,f_{k}$.

For each $i\in\left\{  1,2,\ldots,k\right\}  $, define a map $g_{i}:U\times
V\rightarrow U\times V$ by $g_{i}=f_{i}\times\operatorname*{id}\nolimits_{V}$.

Let $p\in\left\{  1,2,\ldots,k\right\}  $. Thus, $g_{p}=f_{p}\times
\operatorname*{id}\nolimits_{V}$ (by the definition of $g_{p}$). But $f_{p}$
is a transposition of $U$, and thus is a permutation of $U$. Meanwhile,
$\operatorname*{id}\nolimits_{V}$ is a permutation of $V$. Hence, $f_{p}%
\times\operatorname*{id}\nolimits_{V}$ is a permutation of $U\times V$ (by
Corollary \ref{cor.sol.sigmacrosstau.permxperm} (applied to $f_{p}$ and
$\operatorname*{id}\nolimits_{V}$ instead of $\sigma$ and $\tau$)). In other
words, $g_{p}$ is a permutation of $U\times V$ (since $g_{p}=f_{p}%
\times\operatorname*{id}\nolimits_{V}$).

The map $f_{p}$ is a transposition of $U$. In other words, $f_{p}=t_{i,j}$ for
some two distinct elements $i$ and $j$ of $U$ (by the definition of a
transposition of $U$). Consider these $i$ and $j$. Corollary
\ref{cor.sol.sigmacrosstau.transposX.cross3} shows that the permutation
$t_{i,j}\times\operatorname*{id}\nolimits_{V}$ of $U\times V$ satisfies
$\left(  -1\right)  ^{t_{i,j}\times\operatorname*{id}\nolimits_{V}}=\left(
-1\right)  ^{\left\vert V\right\vert }$. Now, $f_{p}=t_{i,j}$. Thus,
$f_{p}\times\operatorname*{id}\nolimits_{V}=t_{i,j}\times\operatorname*{id}%
\nolimits_{V}$. Hence, $g_{p}=f_{p}\times\operatorname*{id}\nolimits_{V}%
=t_{i,j}\times\operatorname*{id}\nolimits_{V}$. Therefore,%
\begin{equation}
\left(  -1\right)  ^{g_{p}}=\left(  -1\right)  ^{t_{i,j}\times
\operatorname*{id}\nolimits_{V}}=\left(  -1\right)  ^{\left\vert V\right\vert
}. \label{pf.cor.sol.sigmacrosstau.sign-sigma-cross-id.3}%
\end{equation}

Now, forget that we fixed $p$. We thus have proven that for each $p\in\left\{
1,2,\ldots,k\right\}  $, the map $g_{p}$ is a permutation of $U\times V$ and
satisfies (\ref{pf.cor.sol.sigmacrosstau.sign-sigma-cross-id.3}).

The maps $g_{1},g_{2},\ldots,g_{k}$ are $k$ permutations of $U\times V$ (since
for each $p\in\left\{  1,2,\ldots,k\right\}  $, the map $g_{p}$ is a
permutation of $U\times V$). Also, the set $U\times V$ is finite (since it is
the product of the finite sets $U$ and $V$). Thus, Proposition
\ref{prop.sol.perm.transX.sign.prod-of-many} (applied to $X=U\times V$ and
$\sigma_{h}=g_{h}$) yields%
\begin{equation}
\left(  -1\right)  ^{g_{1}\circ g_{2}\circ\cdots\circ g_{k}}=\left(
-1\right)  ^{g_{1}}\cdot\left(  -1\right)  ^{g_{2}}\cdot\cdots\cdot\left(
-1\right)  ^{g_{k}}. \label{pf.cor.sol.sigmacrosstau.sign-sigma-cross-id.5}%
\end{equation}

Corollary \ref{cor.sol.sigmacrosstau.fixid} yields%
\[
g_{1}\circ g_{2}\circ\cdots\circ g_{k}=\underbrace{\left(  f_{1}\circ
f_{2}\circ\cdots\circ f_{k}\right)  }_{\substack{=\sigma\\\text{(since }%
\sigma=f_{1}\circ f_{2}\circ\cdots\circ f_{k}\text{)}}}\times
\operatorname*{id}\nolimits_{V}=\sigma\times\operatorname*{id}\nolimits_{V}.
\]
Thus, $\sigma\times\operatorname*{id}\nolimits_{V}=g_{1}\circ g_{2}\circ
\cdots\circ g_{k}$. Therefore,%
\begin{align*}
\left(  -1\right)  ^{\sigma\times\operatorname*{id}\nolimits_{V}}  &  =\left(
-1\right)  ^{g_{1}\circ g_{2}\circ\cdots\circ g_{k}}=\left(  -1\right)
^{g_{1}}\cdot\left(  -1\right)  ^{g_{2}}\cdot\cdots\cdot\left(  -1\right)
^{g_{k}}\\
&  =\prod_{p=1}^{k}\underbrace{\left(  -1\right)  ^{g_{p}}}%
_{\substack{=\left(  -1\right)  ^{\left\vert V\right\vert }\\\text{(by
(\ref{pf.cor.sol.sigmacrosstau.sign-sigma-cross-id.3}))}}}=\prod_{p=1}%
^{k}\left(  -1\right)  ^{\left\vert V\right\vert }=\left(  \left(  -1\right)
^{\left\vert V\right\vert }\right)  ^{k}\\
&  =\left(  -1\right)  ^{\left\vert V\right\vert \cdot k}=\left(  -1\right)
^{k\cdot\left\vert V\right\vert }=\left(  \underbrace{\left(  -1\right)  ^{k}%
}_{\substack{=\left(  -1\right)  ^{\sigma}\\\text{(by
(\ref{pf.cor.sol.sigmacrosstau.sign-sigma-cross-id.1}))}}}\right)
^{\left\vert V\right\vert }=\left(  \left(  -1\right)  ^{\sigma}\right)
^{\left\vert V\right\vert }.
\end{align*}
This proves Corollary \ref{cor.sol.sigmacrosstau.sign-sigma-cross-id}.
\end{proof}
\end{verlong}

\begin{corollary}
\label{cor.sol.sigmacrosstau.sign-id-cross-tau}Let $U$ and $V$ be two finite
sets. Let $\tau$ be a permutation of $V$. Then, the permutation
$\operatorname*{id}\nolimits_{U}\times\tau$ of $U\times V$ satisfies $\left(
-1\right)  ^{\operatorname*{id}\nolimits_{U}\times\tau}=\left(  \left(
-1\right)  ^{\tau}\right)  ^{\left\vert U\right\vert }$.
\end{corollary}

\begin{verlong}
\begin{proof}
[Proof of Corollary \ref{cor.sol.sigmacrosstau.sign-id-cross-tau}.]Corollary
\ref{cor.sol.sigmacrosstau.signs-fggf} (applied to $\sigma=\operatorname*{id}%
\nolimits_{U}$) yields $\left(  -1\right)  ^{\operatorname*{id}\nolimits_{U}%
\times\tau}=\left(  -1\right)  ^{\tau\times\operatorname*{id}\nolimits_{U}}$.
But Corollary \ref{cor.sol.sigmacrosstau.sign-sigma-cross-id} (applied to $V$,
$U$ and $\tau$ instead of $U$, $V$ and $\sigma$) shows that the permutation
$\tau\times\operatorname*{id}\nolimits_{U}$ of $V\times U$ satisfies $\left(
-1\right)  ^{\tau\times\operatorname*{id}\nolimits_{U}}=\left(  \left(
-1\right)  ^{\tau}\right)  ^{\left\vert U\right\vert }$. Hence, $\left(
-1\right)  ^{\operatorname*{id}\nolimits_{U}\times\tau}=\left(  -1\right)
^{\tau\times\operatorname*{id}\nolimits_{U}}=\left(  \left(  -1\right)
^{\tau}\right)  ^{\left\vert U\right\vert }$. This proves Corollary
\ref{cor.sol.sigmacrosstau.sign-id-cross-tau}.
\end{proof}
\end{verlong}

\begin{vershort}
\begin{proof}
[Proof of Corollary \ref{cor.sol.sigmacrosstau.sign-id-cross-tau}.]The map
$\operatorname*{id}\nolimits_{U}$ is a permutation of $U$, whereas $\tau$ is a
permutation of $V$. Hence, $\operatorname*{id}\nolimits_{U}\times\tau$ is a
permutation of $U\times V$ (by Corollary \ref{cor.sol.sigmacrosstau.permxperm}
(applied to $\sigma=\operatorname*{id}\nolimits_{U}$)).

Corollary \ref{cor.sol.sigmacrosstau.signs-fggf} (applied to $\sigma
=\operatorname*{id}\nolimits_{U}$) yields $\left(  -1\right)
^{\operatorname*{id}\nolimits_{U}\times\tau}=\left(  -1\right)  ^{\tau
\times\operatorname*{id}\nolimits_{U}}$. But Corollary
\ref{cor.sol.sigmacrosstau.sign-sigma-cross-id} (applied to $V$, $U$ and
$\tau$ instead of $U$, $V$ and $\sigma$) shows that the permutation
$\tau\times\operatorname*{id}\nolimits_{U}$ of $V\times U$ satisfies $\left(
-1\right)  ^{\tau\times\operatorname*{id}\nolimits_{U}}=\left(  \left(
-1\right)  ^{\tau}\right)  ^{\left\vert U\right\vert }$. Hence, $\left(
-1\right)  ^{\operatorname*{id}\nolimits_{U}\times\tau}=\left(  -1\right)
^{\tau\times\operatorname*{id}\nolimits_{U}}=\left(  \left(  -1\right)
^{\tau}\right)  ^{\left\vert U\right\vert }$. This proves Corollary
\ref{cor.sol.sigmacrosstau.sign-id-cross-tau}.
\end{proof}
\end{vershort}

\begin{corollary}
\label{cor.sol.sigmacrosstau.sign-sigma-cross-tau}Let $U$ and $V$ be two
finite sets. Let $\sigma$ be a permutation of $U$. Let $\tau$ be a permutation
of $V$. Then, the permutation $\sigma\times\tau$ of $U\times V$ satisfies
$\sigma\times\tau=\left(  \sigma\times\operatorname*{id}\nolimits_{V}\right)
\circ\left(  \operatorname*{id}\nolimits_{U}\times\tau\right)  $ and $\left(
-1\right)  ^{\sigma\times\tau}=\left(  \left(  -1\right)  ^{\sigma}\right)
^{\left\vert V\right\vert }\left(  \left(  -1\right)  ^{\tau}\right)
^{\left\vert U\right\vert }$.
\end{corollary}

\begin{vershort}
\begin{proof}
[Proof of Corollary \ref{cor.sol.sigmacrosstau.sign-sigma-cross-tau}%
.]Corollary \ref{cor.sol.sigmacrosstau.permxperm} (applied to $\sigma
=\operatorname*{id}\nolimits_{U}$) shows that $\operatorname*{id}%
\nolimits_{U}\times\tau$ is a permutation of $U\times V$. Similarly,
$\sigma\times\operatorname*{id}\nolimits_{V}$ is a permutation of $U\times V$.

Proposition \ref{prop.sol.sigmacrosstau.fgf'g'} (applied to $U$, $U$, $U$,
$V$, $V$, $V$, $\operatorname*{id}\nolimits_{U}$, $\sigma$, $\tau$ and
$\operatorname*{id}\nolimits_{V}$ instead of $X$, $X^{\prime}$, $X^{\prime
\prime}$, $Y$, $Y^{\prime}$, $Y^{\prime\prime}$, $\alpha$, $\alpha^{\prime}$,
$\beta$ and $\beta^{\prime}$) shows that%
\[
\left(  \sigma\times\operatorname*{id}\nolimits_{V}\right)  \circ\left(
\operatorname*{id}\nolimits_{U}\times\tau\right)  =\underbrace{\left(
\sigma\circ\operatorname*{id}\nolimits_{U}\right)  }_{=\sigma}\times
\underbrace{\left(  \operatorname*{id}\nolimits_{V}\circ\tau\right)  }_{=\tau
}=\sigma\times\tau.
\]
Thus, $\sigma\times\tau=\left(  \sigma\times\operatorname*{id}\nolimits_{V}%
\right)  \circ\left(  \operatorname*{id}\nolimits_{U}\times\tau\right)  $.
Hence,%
\begin{align*}
\left(  -1\right)  ^{\sigma\times\tau}  &  =\left(  -1\right)  ^{\left(
\sigma\times\operatorname*{id}\nolimits_{V}\right)  \circ\left(
\operatorname*{id}\nolimits_{U}\times\tau\right)  }=\underbrace{\left(
-1\right)  ^{\sigma\times\operatorname*{id}\nolimits_{V}}}_{\substack{=\left(
\left(  -1\right)  ^{\sigma}\right)  ^{\left\vert V\right\vert }\\\text{(by
Corollary \ref{cor.sol.sigmacrosstau.sign-sigma-cross-id})}}}\cdot
\underbrace{\left(  -1\right)  ^{\operatorname*{id}\nolimits_{U}\times\tau}%
}_{\substack{=\left(  \left(  -1\right)  ^{\tau}\right)  ^{\left\vert
U\right\vert }\\\text{(by Corollary
\ref{cor.sol.sigmacrosstau.sign-id-cross-tau})}}}\\
&  \ \ \ \ \ \ \ \ \ \ \left(
\begin{array}
[c]{c}%
\text{by Exercise \ref{exe.ps4.2} \textbf{(c)} (applied to }U\times V\text{,
}\sigma\times\operatorname*{id}\nolimits_{V}\\
\text{and }\operatorname*{id}\nolimits_{U}\times\tau\text{ instead of
}X\text{, }\sigma\text{ and }\tau\text{)}%
\end{array}
\right) \\
&  =\left(  \left(  -1\right)  ^{\sigma}\right)  ^{\left\vert V\right\vert
}\left(  \left(  -1\right)  ^{\tau}\right)  ^{\left\vert U\right\vert }.
\end{align*}
This completes the proof of Corollary
\ref{cor.sol.sigmacrosstau.sign-sigma-cross-tau}.
\end{proof}
\end{vershort}

\begin{verlong}
\begin{proof}
[Proof of Corollary \ref{cor.sol.sigmacrosstau.sign-sigma-cross-tau}.]The set
$U\times V$ is finite (since it is the product of the two finite sets $U$ and
$V$). The map $\sigma\times\tau$ is a permutation of $U\times V$ (by Corollary
\ref{cor.sol.sigmacrosstau.permxperm}).

The map $\operatorname*{id}\nolimits_{U}$ is a permutation of $U$, whereas
$\tau$ is a permutation of $V$. Hence, $\operatorname*{id}\nolimits_{U}%
\times\tau$ is a permutation of $U\times V$ (by Corollary
\ref{cor.sol.sigmacrosstau.permxperm} (applied to $\sigma=\operatorname*{id}%
\nolimits_{U}$)).

The map $\sigma$ is a permutation of $U$, whereas $\operatorname*{id}%
\nolimits_{V}$ is a permutation of $V$. Hence, $\sigma\times\operatorname*{id}%
\nolimits_{V}$ is a permutation of $U\times V$ (by Corollary
\ref{cor.sol.sigmacrosstau.permxperm} (applied to $\tau=\operatorname*{id}%
\nolimits_{V}$)).

Proposition \ref{prop.sol.sigmacrosstau.fgf'g'} (applied to $U$, $U$, $U$,
$V$, $V$, $V$, $\operatorname*{id}\nolimits_{U}$, $\sigma$, $\tau$ and
$\operatorname*{id}\nolimits_{V}$ instead of $X$, $X^{\prime}$, $X^{\prime
\prime}$, $Y$, $Y^{\prime}$, $Y^{\prime\prime}$, $\alpha$, $\alpha^{\prime}$,
$\beta$ and $\beta^{\prime}$) shows that%
\[
\left(  \sigma\times\operatorname*{id}\nolimits_{V}\right)  \circ\left(
\operatorname*{id}\nolimits_{U}\times\tau\right)  =\underbrace{\left(
\sigma\circ\operatorname*{id}\nolimits_{U}\right)  }_{=\sigma}\times
\underbrace{\left(  \operatorname*{id}\nolimits_{V}\circ\tau\right)  }_{=\tau
}=\sigma\times\tau.
\]
Thus, $\sigma\times\tau=\left(  \sigma\times\operatorname*{id}\nolimits_{V}%
\right)  \circ\left(  \operatorname*{id}\nolimits_{U}\times\tau\right)  $.
Hence,%
\begin{align*}
\left(  -1\right)  ^{\sigma\times\tau}  &  =\left(  -1\right)  ^{\left(
\sigma\times\operatorname*{id}\nolimits_{V}\right)  \circ\left(
\operatorname*{id}\nolimits_{U}\times\tau\right)  }=\underbrace{\left(
-1\right)  ^{\sigma\times\operatorname*{id}\nolimits_{V}}}_{\substack{=\left(
\left(  -1\right)  ^{\sigma}\right)  ^{\left\vert V\right\vert }\\\text{(by
Corollary \ref{cor.sol.sigmacrosstau.sign-sigma-cross-id})}}}\cdot
\underbrace{\left(  -1\right)  ^{\operatorname*{id}\nolimits_{U}\times\tau}%
}_{\substack{=\left(  \left(  -1\right)  ^{\tau}\right)  ^{\left\vert
U\right\vert }\\\text{(by Corollary
\ref{cor.sol.sigmacrosstau.sign-id-cross-tau})}}}\\
&  \ \ \ \ \ \ \ \ \ \ \left(
\begin{array}
[c]{c}%
\text{by Exercise \ref{exe.ps4.2} \textbf{(c)} (applied to }U\times V\text{,
}\sigma\times\operatorname*{id}\nolimits_{V}\\
\text{and }\operatorname*{id}\nolimits_{U}\times\tau\text{ instead of
}X\text{, }\sigma\text{ and }\tau\text{)}%
\end{array}
\right) \\
&  =\left(  \left(  -1\right)  ^{\sigma}\right)  ^{\left\vert V\right\vert
}\left(  \left(  -1\right)  ^{\tau}\right)  ^{\left\vert U\right\vert }.
\end{align*}
This completes the proof of Corollary
\ref{cor.sol.sigmacrosstau.sign-sigma-cross-tau}.
\end{proof}
\end{verlong}

\begin{vershort}
\begin{proof}
[Solution to Exercise \ref{exe.perm.sigmacrosstau}.]Exercise
\ref{exe.perm.sigmacrosstau} \textbf{(a)} follows from Corollary
\ref{cor.sol.sigmacrosstau.permxperm}. Parts \textbf{(b)} and \textbf{(c)} of
Exercise \ref{exe.perm.sigmacrosstau} follow from Corollary
\ref{cor.sol.sigmacrosstau.sign-sigma-cross-tau}.
\end{proof}
\end{vershort}

\begin{verlong}
\begin{proof}
[Solution to Exercise \ref{exe.perm.sigmacrosstau}.]\textbf{(a)} The map
$\sigma\times\tau$ is a permutation of $U\times V$ (by Corollary
\ref{cor.sol.sigmacrosstau.permxperm}). This solves Exercise
\ref{exe.perm.sigmacrosstau} \textbf{(a)}.

Corollary \ref{cor.sol.sigmacrosstau.sign-sigma-cross-tau} yields
$\sigma\times\tau=\left(  \sigma\times\operatorname*{id}\nolimits_{V}\right)
\circ\left(  \operatorname*{id}\nolimits_{U}\times\tau\right)  $ and
\[
\left(  -1\right)  ^{\sigma\times\tau}=\left(  \left(  -1\right)  ^{\sigma
}\right)  ^{\left\vert V\right\vert }\left(  \left(  -1\right)  ^{\tau
}\right)  ^{\left\vert U\right\vert } .
\]

We have $\sigma\times\tau=\left(  \sigma\times\operatorname*{id}%
\nolimits_{V}\right)  \circ\left(  \operatorname*{id}\nolimits_{U}\times
\tau\right)  $; this solves Exercise \ref{exe.perm.sigmacrosstau} \textbf{(b)}.

We have $\left(  -1\right)  ^{\sigma\times\tau}=\left(  \left(  -1\right)
^{\sigma}\right)  ^{\left\vert V\right\vert }\left(  \left(  -1\right)
^{\tau}\right)  ^{\left\vert U\right\vert }$. This solves Exercise
\ref{exe.perm.sigmacrosstau} \textbf{(c)}.
\end{proof}
\end{verlong}

\subsection{\label{sect.sol.perm.footrule}Solution to Exercise
\ref{exe.perm.footrule}}

Throughout Section \ref{sect.sol.perm.footrule}, we shall use the same
notations that were used in Section \ref{sect.perm.lehmer}. We shall also use
Definition \ref{def.iverson}.

We begin with several minor results that we shall use in our solution of
Exercise \ref{exe.perm.footrule}. The first of these is a well-known
inequality, known as the \textit{triangle inequality for numbers}:

\begin{theorem}
\label{thm.ineq.triangle-R}Let $x$ and $y$ be two rational numbers (or real
numbers). Then, $\left\vert x\right\vert +\left\vert y\right\vert
\geq\left\vert x+y\right\vert $.
\end{theorem}

\begin{proof}
[Proof of Theorem \ref{thm.ineq.triangle-R}.]Every rational number (or real
number) $z$ satisfies%
\begin{equation}
\left\vert z\right\vert =\max\left\{  z,-z\right\}  \geq z
\label{pf.thm.ineq.triangle-R.1}%
\end{equation}
and%
\begin{equation}
\left\vert z\right\vert =\max\left\{  z,-z\right\}  \geq-z.
\label{pf.thm.ineq.triangle-R.2}%
\end{equation}

Now, $x+y$ is a rational (or real) number, and thus is either $\geq0$ or $<0$.
Hence, we are in one of the following two cases:

\textit{Case 1:} We have $x+y\geq0$.

\textit{Case 2:} We have $x+y<0$.

\begin{vershort}
Let us first consider Case 1. In this case, we have $x+y\geq0$. Hence,
$\left\vert x+y\right\vert =x+y$. Now, $\underbrace{\left\vert x\right\vert
}_{\substack{\geq x\\\text{(by (\ref{pf.thm.ineq.triangle-R.1}))}%
}}+\underbrace{\left\vert y\right\vert }_{\substack{\geq y\\\text{(by
(\ref{pf.thm.ineq.triangle-R.1}))}}}\geq x+y=\left\vert x+y\right\vert $.
Hence, Theorem \ref{thm.ineq.triangle-R} is proven in Case 1.

Let us now consider Case 2. In this case, we have $x+y<0$. Hence, $\left\vert
x+y\right\vert =-\left(  x+y\right)  =\left(  -x\right)  +\left(  -y\right)
$. Now, $\underbrace{\left\vert x\right\vert }_{\substack{\geq-x\\\text{(by
(\ref{pf.thm.ineq.triangle-R.2}))}}}+\underbrace{\left\vert y\right\vert
}_{\substack{\geq-y\\\text{(by (\ref{pf.thm.ineq.triangle-R.2}))}}}\geq\left(
-x\right)  +\left(  -y\right)  =\left\vert x+y\right\vert $. Hence, Theorem
\ref{thm.ineq.triangle-R} is proven in Case 2.
\end{vershort}

\begin{verlong}
Let us first consider Case 1. In this case, we have $x+y\geq0$. Hence,
$\left\vert x+y\right\vert =x+y$. Thus, $x+y=\left\vert x+y\right\vert $.
Now,
\[
\underbrace{\left\vert x\right\vert }_{\substack{\geq x\\\text{(by
(\ref{pf.thm.ineq.triangle-R.1}) (applied to }z=x\text{))}}%
}+\underbrace{\left\vert y\right\vert }_{\substack{\geq y\\\text{(by
(\ref{pf.thm.ineq.triangle-R.1}) (applied to }z=y\text{))}}}\geq
x+y=\left\vert x+y\right\vert .
\]
Hence, Theorem \ref{thm.ineq.triangle-R} is proven in Case 1.

Let us now consider Case 2. In this case, we have $x+y<0$. Hence, $\left\vert
x+y\right\vert =-\left(  x+y\right)  =\left(  -x\right)  +\left(  -y\right)
$. Thus, $\left(  -x\right)  +\left(  -y\right)  =\left\vert x+y\right\vert $.
Now,
\[
\underbrace{\left\vert x\right\vert }_{\substack{\geq-x\\\text{(by
(\ref{pf.thm.ineq.triangle-R.2}) (applied to }z=x\text{))}}%
}+\underbrace{\left\vert y\right\vert }_{\substack{\geq-y\\\text{(by
(\ref{pf.thm.ineq.triangle-R.2}) (applied to }z=y\text{))}}}\geq\left(
-x\right)  +\left(  -y\right)  =\left\vert x+y\right\vert .
\]
Hence, Theorem \ref{thm.ineq.triangle-R} is proven in Case 2.
\end{verlong}

\begin{vershort}
We have now proven Theorem \ref{thm.ineq.triangle-R} in each of the two Cases
1 and 2. Hence, Theorem \ref{thm.ineq.triangle-R} always holds. \qedhere

\end{vershort}

\begin{verlong}
We have now proven Theorem \ref{thm.ineq.triangle-R} in each of the two Cases
1 and 2. Since these two Cases cover all possibilities, we thus conclude that
Theorem \ref{thm.ineq.triangle-R} always holds.
\end{verlong}
\end{proof}

Another simple property of summations will be useful:

\begin{lemma}
\label{lem.sol.perm.footrule.1}Let $n\in\mathbb{N}$. Let $u$ and $v$ be two
elements of $\left[  n\right]  $.

\textbf{(a)} We have $1=\left[  u<v\right]  +\left[  v\leq u\right]  $.

\textbf{(b)} We have $\sum_{k\in\left[  n\right]  }\left[  k>u\right]  \left[
k<v\right]  =\left[  u<v\right]  \left(  v-u-1\right)  $.

\textbf{(c)} We have $\sum_{k\in\left[  n\right]  }\left[  k>u\right]  \left[
k\leq v\right]  =\left[  u<v\right]  \left(  v-u\right)  $.
\end{lemma}

\begin{proof}
[Proof of Lemma \ref{lem.sol.perm.footrule.1}.]\textbf{(a)} The logical
statements $\left(  v\leq u\right)  $ and $\left(  \text{not }u<v\right)  $
are equivalent (since $\left(  v\leq u\right)  \Longleftrightarrow\left(
u\geq v\right)  \Longleftrightarrow\left(  \text{not }u<v\right)  $). Hence,
Exercise \ref{exe.iverson-prop} \textbf{(a)} (applied to $\mathcal{A}=\left(
v\leq u\right)  $ and $\mathcal{B}=\left(  \text{not }u<v\right)  $) shows
that $\left[  v\leq u\right]  =\left[  \text{not }u<v\right]  =1-\left[
u<v\right]  $ (by Exercise \ref{exe.iverson-prop} \textbf{(b)} (applied to
$\mathcal{A}=\left(  u<v\right)  $)). Adding $\left[  u<v\right]  $ to both
sides of this equality, we obtain $\left[  u<v\right]  +\left[  v\leq
u\right]  =1$. This proves Lemma \ref{lem.sol.perm.footrule.1} \textbf{(a)}.

\begin{vershort}
\textbf{(b)} We have $u\in\left[  n\right]  $, so that $1\leq u\leq n$.
Similarly, $1\leq v\leq n$.
\end{vershort}

\begin{verlong}
\textbf{(b)} We have $u\in\left[  n\right]  =\left\{  1,2,\ldots,n\right\}  $
(by the definition of $\left[  n\right]  $), so that $1\leq u\leq n$. The same
argument (applied to $v$ instead of $u$) yields $1\leq v\leq n$.
\end{verlong}

\begin{verlong}
We have $\left[  n\right]  =\left\{  1,2,\ldots,n\right\}  $ and thus
$\sum_{k\in\left[  n\right]  }=\sum_{k\in\left\{  1,2,\ldots,n\right\}  }%
=\sum_{k=1}^{n}$ (an equality between summation signs).
\end{verlong}

Now,%
\begin{align}
&  \underbrace{\sum_{k\in\left[  n\right]  }}_{=\sum_{k=1}^{n}}\left[
k>u\right]  \left[  k<v\right] \nonumber\\
&  =\sum_{k=1}^{n}\left[  k>u\right]  \left[  k<v\right]  =\sum_{k=1}%
^{u}\underbrace{\left[  k>u\right]  }_{\substack{=0\\\text{(since we don't
have }k>u\\\text{(because }k\leq u\text{))}}}\left[  k<v\right]  +\sum
_{k=u+1}^{n}\underbrace{\left[  k>u\right]  }_{\substack{=1\\\text{(since
}k\geq u+1>u\text{)}}}\left[  k<v\right] \nonumber\\
&  \ \ \ \ \ \ \ \ \ \ \left(  \text{here, we have split the sum at
}k=u\text{, since }1\leq u\leq n\right) \nonumber\\
&  =\underbrace{\sum_{k=1}^{u}0\left[  k<v\right]  }_{=0}+\sum_{k=u+1}%
^{n}\left[  k<v\right]  =\sum_{k=u+1}^{n}\left[  k<v\right]  .
\label{pf.lem.sol.perm.footrule.1.b.1}%
\end{align}

We are in one of the following two cases:

\textit{Case 1:} We have $u<v$.

\textit{Case 2:} We don't have $u<v$.

Let us first consider Case 1. In this case, we have $u<v$. Thus, $u\leq v-1$
(since $u$ and $v$ are integers), so that $u+1\leq v$. Now,
(\ref{pf.lem.sol.perm.footrule.1.b.1}) becomes%
\begin{align*}
&  \sum_{k\in\left[  n\right]  }\left[  k>u\right]  \left[  k<v\right] \\
&  =\sum_{k=u+1}^{n}\left[  k<v\right]  =\sum_{k=u+1}^{v-1}\underbrace{\left[
k<v\right]  }_{\substack{=1\\\text{(since }k\leq v-1<v\text{)}}}+\sum
_{k=v}^{n}\underbrace{\left[  k<v\right]  }_{\substack{=0\\\text{(since we
don't have }k<v\\\text{(because }k\geq v\text{))}}}\\
&  \ \ \ \ \ \ \ \ \ \ \left(  \text{here, we have split the sum at
}k=v-1\text{, since }u+1\leq v\leq n\right) \\
&  =\sum_{k=u+1}^{v-1}1+\underbrace{\sum_{k=v}^{n}0}_{=0}=\sum_{k=u+1}%
^{v-1}1=\left(  \left(  v-1\right)  -u\right)  \cdot
1\ \ \ \ \ \ \ \ \ \ \left(  \text{since }u\leq v-1\right) \\
&  =v-u-1.
\end{align*}
Comparing this with%
\[
\underbrace{\left[  u<v\right]  }_{\substack{=1\\\text{(since }u<v\text{)}%
}}\left(  v-u-1\right)  =v-u-1,
\]
we obtain $\sum_{k\in\left[  n\right]  }\left[  k>u\right]  \left[
k<v\right]  =\left[  u<v\right]  \left(  v-u-1\right)  $. Thus, Lemma
\ref{lem.sol.perm.footrule.1} \textbf{(b)} is proven in Case 1.

Let us next consider Case 2. In this case, we don't have $u<v$. Thus, we have
$u\geq v$. Now, (\ref{pf.lem.sol.perm.footrule.1.b.1}) becomes%
\[
\sum_{k\in\left[  n\right]  }\left[  k>u\right]  \left[  k<v\right]
=\sum_{k=u+1}^{n}\underbrace{\left[  k<v\right]  }_{\substack{=0\\\text{(since
we don't have }k<v\\\text{(because }k\geq u+1\geq u\geq v\text{))}}%
}=\sum_{k=u+1}^{n}0=0.
\]
Comparing this with%
\[
\underbrace{\left[  u<v\right]  }_{\substack{=0\\\text{(since we don't have
}u<v\text{)}}}\left(  v-u-1\right)  =0\left(  v-u-1\right)  =0,
\]
we obtain $\sum_{k\in\left[  n\right]  }\left[  k>u\right]  \left[
k<v\right]  =\left[  u<v\right]  \left(  v-u-1\right)  $. Thus, Lemma
\ref{lem.sol.perm.footrule.1} \textbf{(b)} is proven in Case 2.

\begin{vershort}
We have now proven Lemma \ref{lem.sol.perm.footrule.1} \textbf{(b)} in each of
the two Cases 1 and 2. Hence, Lemma \ref{lem.sol.perm.footrule.1} \textbf{(b)}
always holds.
\end{vershort}

\begin{verlong}
We have now proven Lemma \ref{lem.sol.perm.footrule.1} \textbf{(b)} in each of
the two Cases 1 and 2. Since these two Cases cover all possibilities, we thus
conclude that Lemma \ref{lem.sol.perm.footrule.1} \textbf{(b)} always holds.
\end{verlong}

\begin{vershort}
\textbf{(c)} We have $u\in\left[  n\right]  $, so that $1\leq u\leq n$.
Similarly, $1\leq v\leq n$.
\end{vershort}

\begin{verlong}
\textbf{(c)} We have $u\in\left[  n\right]  =\left\{  1,2,\ldots,n\right\}  $
(by the definition of $\left[  n\right]  $), so that $1\leq u\leq n$. The same
argument (applied to $v$ instead of $u$) yields $1\leq v\leq n$.
\end{verlong}

\begin{verlong}
We have $\left[  n\right]  =\left\{  1,2,\ldots,n\right\}  $ and thus
$\sum_{k\in\left[  n\right]  }=\sum_{k\in\left\{  1,2,\ldots,n\right\}  }%
=\sum_{k=1}^{n}$ (an equality between summation signs).
\end{verlong}

\begin{vershort}
Now,%
\begin{equation}
\sum_{k\in\left[  n\right]  }\left[  k>u\right]  \left[  k\leq v\right]
=\sum_{k=u+1}^{n}\left[  k\leq v\right]  .
\label{pf.lem.sol.perm.footrule.1.c.short.1}%
\end{equation}
(This can be proven in the exact same way as
(\ref{pf.lem.sol.perm.footrule.1.b.1}), except that all the \textquotedblleft%
$\left[  k<v\right]  $\textquotedblright\ terms must be replaced by
\textquotedblleft$\left[  k\leq v\right]  $\textquotedblright.)
\end{vershort}

\begin{verlong}
Now,%
\begin{align}
&  \underbrace{\sum_{k\in\left[  n\right]  }}_{=\sum_{k=1}^{n}}\left[
k>u\right]  \left[  k\leq v\right] \nonumber\\
&  =\sum_{k=1}^{n}\left[  k>u\right]  \left[  k\leq v\right]  =\sum_{k=1}%
^{u}\underbrace{\left[  k>u\right]  }_{\substack{=0\\\text{(since we don't
have }k>u\\\text{(because }k\leq u\text{))}}}\left[  k\leq v\right]
+\sum_{k=u+1}^{n}\underbrace{\left[  k>u\right]  }_{\substack{=1\\\text{(since
}k\geq u+1>u\text{)}}}\left[  k\leq v\right] \nonumber\\
&  \ \ \ \ \ \ \ \ \ \ \left(  \text{here, we have split the sum at
}k=u\text{, since }1\leq u\leq n\right) \nonumber\\
&  =\underbrace{\sum_{k=1}^{u}0\left[  k\leq v\right]  }_{=0}+\sum_{k=u+1}%
^{n}\left[  k\leq v\right]  =\sum_{k=u+1}^{n}\left[  k\leq v\right]  .
\label{pf.lem.sol.perm.footrule.1.c.1}%
\end{align}

\end{verlong}

We are in one of the following two cases:

\textit{Case 1:} We have $u<v$.

\textit{Case 2:} We don't have $u<v$.

\begin{vershort}
Let us first consider Case 1. In this case, we have $u<v$. Thus, $u\leq v-1$
(since $u$ and $v$ are integers), so that $u+1\leq v$. Now,
(\ref{pf.lem.sol.perm.footrule.1.c.short.1}) becomes%
\begin{align*}
&  \sum_{k\in\left[  n\right]  }\left[  k>u\right]  \left[  k\leq v\right] \\
&  =\sum_{k=u+1}^{n}\left[  k\leq v\right]  =\sum_{k=u+1}^{v}%
\underbrace{\left[  k\leq v\right]  }_{\substack{=1\\\text{(since }k\leq
v\text{)}}}+\sum_{k=v+1}^{n}\underbrace{\left[  k\leq v\right]  }%
_{\substack{=0\\\text{(since we don't have }k\leq v\\\text{(because }k\geq
v+1>v\text{))}}}\\
&  \ \ \ \ \ \ \ \ \ \ \left(  \text{here, we have split the sum at
}k=v\text{, since }u+1\leq v\leq n\right) \\
&  =\sum_{k=u+1}^{v}1+\underbrace{\sum_{k=v+1}^{n}0}_{=0}=\sum_{k=u+1}%
^{v}1=\left(  v-u\right)  \cdot1\ \ \ \ \ \ \ \ \ \ \left(  \text{since }u\leq
v-1\leq v\right) \\
&  =v-u.
\end{align*}
Comparing this with $\underbrace{\left[  u<v\right]  }%
_{\substack{=1\\\text{(since }u<v\text{)}}}\left(  v-u\right)  =v-u$, we
obtain $\sum_{k\in\left[  n\right]  }\left[  k>u\right]  \left[  k\leq
v\right]  =\left[  u<v\right]  \left(  v-u\right)  $. Thus, Lemma
\ref{lem.sol.perm.footrule.1} \textbf{(c)} is proven in Case 1.
\end{vershort}

\begin{verlong}
Let us first consider Case 1. In this case, we have $u<v$. Thus, $u\leq v-1$
(since $u$ and $v$ are integers), so that $u+1\leq v$. Now,
(\ref{pf.lem.sol.perm.footrule.1.c.1}) becomes%
\begin{align*}
&  \sum_{k\in\left[  n\right]  }\left[  k>u\right]  \left[  k\leq v\right] \\
&  =\sum_{k=u+1}^{n}\left[  k\leq v\right]  =\sum_{k=u+1}^{v}%
\underbrace{\left[  k\leq v\right]  }_{\substack{=1\\\text{(since }k\leq
v\text{)}}}+\sum_{k=v+1}^{n}\underbrace{\left[  k\leq v\right]  }%
_{\substack{=0\\\text{(since we don't have }k\leq v\\\text{(because }k\geq
v+1>v\text{))}}}\\
&  \ \ \ \ \ \ \ \ \ \ \left(  \text{here, we have split the sum at
}k=v\text{, since }u+1\leq v\leq n\right) \\
&  =\sum_{k=u+1}^{v}1+\underbrace{\sum_{k=v+1}^{n}0}_{=0}=\sum_{k=u+1}%
^{v}1=\left(  v-u\right)  \cdot1\ \ \ \ \ \ \ \ \ \ \left(  \text{since }u\leq
v-1\leq v\right) \\
&  =v-u.
\end{align*}
Comparing this with%
\[
\underbrace{\left[  u<v\right]  }_{\substack{=1\\\text{(since }u<v\text{)}%
}}\left(  v-u\right)  =v-u,
\]
we obtain $\sum_{k\in\left[  n\right]  }\left[  k>u\right]  \left[  k\leq
v\right]  =\left[  u<v\right]  \left(  v-u\right)  $. Thus, Lemma
\ref{lem.sol.perm.footrule.1} \textbf{(c)} is proven in Case 1.
\end{verlong}

\begin{vershort}
Let us next consider Case 2. In this case, we don't have $u<v$. Thus, we have
$u\geq v$. Now, (\ref{pf.lem.sol.perm.footrule.1.c.short.1}) becomes%
\[
\sum_{k\in\left[  n\right]  }\left[  k>u\right]  \left[  k\leq v\right]
=\sum_{k=u+1}^{n}\underbrace{\left[  k\leq v\right]  }%
_{\substack{=0\\\text{(since we don't have }k\leq v\\\text{(because }k\geq
u+1>u\geq v\text{))}}}=\sum_{k=u+1}^{n}0=0.
\]
Comparing this with $\underbrace{\left[  u<v\right]  }%
_{\substack{=0\\\text{(since we don't have }u<v\text{)}}}\left(  v-u\right)
=0\left(  v-u\right)  =0$, we obtain \newline$\sum_{k\in\left[  n\right]
}\left[  k>u\right]  \left[  k\leq v\right]  =\left[  u<v\right]  \left(
v-u\right)  $. Thus, Lemma \ref{lem.sol.perm.footrule.1} \textbf{(c)} is
proven in Case 2.
\end{vershort}

\begin{verlong}
Let us next consider Case 2. In this case, we don't have $u<v$. Thus, we have
$u\geq v$. Now, (\ref{pf.lem.sol.perm.footrule.1.c.1}) becomes%
\[
\sum_{k\in\left[  n\right]  }\left[  k>u\right]  \left[  k\leq v\right]
=\sum_{k=u+1}^{n}\underbrace{\left[  k\leq v\right]  }%
_{\substack{=0\\\text{(since we don't have }k\leq v\\\text{(because }k\geq
u+1>u\geq v\text{))}}}=\sum_{k=u+1}^{n}0=0.
\]
Comparing this with%
\[
\underbrace{\left[  u<v\right]  }_{\substack{=0\\\text{(since we don't have
}u<v\text{)}}}\left(  v-u\right)  =0\left(  v-u\right)  =0,
\]
we obtain $\sum_{k\in\left[  n\right]  }\left[  k>u\right]  \left[  k\leq
v\right]  =\left[  u<v\right]  \left(  v-u\right)  $. Thus, Lemma
\ref{lem.sol.perm.footrule.1} \textbf{(c)} is proven in Case 2.
\end{verlong}

\begin{vershort}
We have now proven Lemma \ref{lem.sol.perm.footrule.1} \textbf{(c)} in each of
the two Cases 1 and 2. Hence, Lemma \ref{lem.sol.perm.footrule.1} \textbf{(c)}
always holds. \qedhere

\end{vershort}

\begin{verlong}
We have now proven Lemma \ref{lem.sol.perm.footrule.1} \textbf{(c)} in each of
the two Cases 1 and 2. Since these two Cases cover all possibilities, we thus
conclude that Lemma \ref{lem.sol.perm.footrule.1} \textbf{(c)} always holds.
\end{verlong}
\end{proof}

\begin{proof}
[Solution to Exercise \ref{exe.perm.footrule}.]We start with some simple observations.

\begin{vershort}
Recall that $S_{n}$ is the set of all permutations of $\left\{  1,2,\ldots
,n\right\}  $. In other words, $S_{n}$ is the set of all permutations of
$\left[  n\right]  $ (since $\left\{  1,2,\ldots,n\right\}  =\left[  n\right]
$). Hence, $\sigma$ is a permutation of $\left[  n\right]  $ (since $\sigma\in
S_{n}$). In other words, $\sigma$ is a bijective map $\left[  n\right]
\rightarrow\left[  n\right]  $.
\end{vershort}

\begin{verlong}
We have $\left[  n\right]  =\left\{  1,2,\ldots,n\right\}  $ (by the
definition of $\left[  n\right]  $).

Recall that $S_{n}$ is the set of all permutations of the set $\left\{
1,2,\ldots,n\right\}  $. In other words, $S_{n}$ is the set of all
permutations of the set $\left[  n\right]  $ (since $\left\{  1,2,\ldots
,n\right\}  =\left[  n\right]  $).

We have $\sigma\in S_{n}$. In other words, $\sigma$ is a permutation of the
set $\left[  n\right]  $ (since $S_{n}$ is the set of all permutations of the
set $\left[  n\right]  $). In other words, $\sigma$ is a bijection $\left[
n\right]  \rightarrow\left[  n\right]  $. In other words, $\sigma$ is a
bijective map $\left[  n\right]  \rightarrow\left[  n\right]  $.
\end{verlong}

Each $i\in\left[  n\right]  $ satisfies exactly one of the three statements
$\sigma\left(  i\right)  <i$, $\sigma\left(  i\right)  =i$ and $\sigma\left(
i\right)  >i$. Hence, we can subdivide the sum $\sum_{i\in\left[  n\right]
}\left(  \sigma\left(  i\right)  -i\right)  $ as follows:%
\begin{align*}
\sum_{i\in\left[  n\right]  }\left(  \sigma\left(  i\right)  -i\right)   &
=\sum_{\substack{i\in\left[  n\right]  ;\\\sigma\left(  i\right)
<i}}\underbrace{\left(  \sigma\left(  i\right)  -i\right)  }_{=-\left(
i-\sigma\left(  i\right)  \right)  }+\sum_{\substack{i\in\left[  n\right]
;\\\sigma\left(  i\right)  =i}}\underbrace{\left(  \sigma\left(  i\right)
-i\right)  }_{\substack{=0\\\text{(since }\sigma\left(  i\right)  =i\text{)}%
}}+\sum_{\substack{i\in\left[  n\right]  ;\\\sigma\left(  i\right)
>i}}\left(  \sigma\left(  i\right)  -i\right) \\
&  =\underbrace{\sum_{\substack{i\in\left[  n\right]  ;\\\sigma\left(
i\right)  <i}}\left(  -\left(  i-\sigma\left(  i\right)  \right)  \right)
}_{=-\sum_{\substack{i\in\left[  n\right]  ;\\\sigma\left(  i\right)
<i}}\left(  i-\sigma\left(  i\right)  \right)  }+\underbrace{\sum
_{\substack{i\in\left[  n\right]  ;\\\sigma\left(  i\right)  =i}}0}_{=0}%
+\sum_{\substack{i\in\left[  n\right]  ;\\\sigma\left(  i\right)  >i}}\left(
\sigma\left(  i\right)  -i\right) \\
&  =-\sum_{\substack{i\in\left[  n\right]  ;\\\sigma\left(  i\right)
<i}}\left(  i-\sigma\left(  i\right)  \right)  +\sum_{\substack{i\in\left[
n\right]  ;\\\sigma\left(  i\right)  >i}}\left(  \sigma\left(  i\right)
-i\right) \\
&  =\sum_{\substack{i\in\left[  n\right]  ;\\\sigma\left(  i\right)
>i}}\left(  \sigma\left(  i\right)  -i\right)  -\sum_{\substack{i\in\left[
n\right]  ;\\\sigma\left(  i\right)  <i}}\left(  i-\sigma\left(  i\right)
\right)  .
\end{align*}
Hence,%
\begin{align*}
&  \sum_{\substack{i\in\left[  n\right]  ;\\\sigma\left(  i\right)
>i}}\left(  \sigma\left(  i\right)  -i\right)  -\sum_{\substack{i\in\left[
n\right]  ;\\\sigma\left(  i\right)  <i}}\left(  i-\sigma\left(  i\right)
\right) \\
&  =\sum_{i\in\left[  n\right]  }\left(  \sigma\left(  i\right)  -i\right)
=\underbrace{\sum_{i\in\left[  n\right]  }\sigma\left(  i\right)
}_{\substack{=\sum_{i\in\left[  n\right]  }i\\\text{(here, we have substituted
}i\\\text{for }\sigma\left(  i\right)  \text{ in the sum, since }%
\sigma\\\text{is a bijection }\left[  n\right]  \rightarrow\left[  n\right]
\text{)}}}-\sum_{i\in\left[  n\right]  }i=\sum_{i\in\left[  n\right]  }%
i-\sum_{i\in\left[  n\right]  }i=0.
\end{align*}
In other words,%
\begin{equation}
\sum_{\substack{i\in\left[  n\right]  ;\\\sigma\left(  i\right)  >i}}\left(
\sigma\left(  i\right)  -i\right)  =\sum_{\substack{i\in\left[  n\right]
;\\\sigma\left(  i\right)  <i}}\left(  i-\sigma\left(  i\right)  \right)  .
\label{sol.exe.perm.footrule.a.1}%
\end{equation}

\textbf{(a)} Each $i\in\left[  n\right]  $ satisfies exactly one of the three
statements $\sigma\left(  i\right)  <i$, $\sigma\left(  i\right)  =i$ and
$\sigma\left(  i\right)  >i$. Hence, we can subdivide the sum $\sum
_{i\in\left[  n\right]  }\left\vert \sigma\left(  i\right)  -i\right\vert $ as
follows:%
\begin{align*}
\sum_{i\in\left[  n\right]  }\left\vert \sigma\left(  i\right)  -i\right\vert
&  =\sum_{\substack{i\in\left[  n\right]  ;\\\sigma\left(  i\right)
<i}}\underbrace{\left\vert \sigma\left(  i\right)  -i\right\vert
}_{\substack{=-\left(  \sigma\left(  i\right)  -i\right)  \\\text{(since
}\sigma\left(  i\right)  -i<0\\\text{(because }\sigma\left(  i\right)
<i\text{))}}}+\sum_{\substack{i\in\left[  n\right]  ;\\\sigma\left(  i\right)
=i}}\left\vert \underbrace{\sigma\left(  i\right)  -i}%
_{\substack{=0\\\text{(since }\sigma\left(  i\right)  =i\text{)}}}\right\vert
+\sum_{\substack{i\in\left[  n\right]  ;\\\sigma\left(  i\right)
>i}}\underbrace{\left\vert \sigma\left(  i\right)  -i\right\vert
}_{\substack{=\sigma\left(  i\right)  -i\\\text{(since }\sigma\left(
i\right)  -i>0\\\text{(because }\sigma\left(  i\right)  >i\text{))}}}\\
&  =\sum_{\substack{i\in\left[  n\right]  ;\\\sigma\left(  i\right)
<i}}\underbrace{\left(  -\left(  \sigma\left(  i\right)  -i\right)  \right)
}_{=i-\sigma\left(  i\right)  }+\sum_{\substack{i\in\left[  n\right]
;\\\sigma\left(  i\right)  =i}}\underbrace{\left\vert 0\right\vert }_{=0}%
+\sum_{\substack{i\in\left[  n\right]  ;\\\sigma\left(  i\right)  >i}}\left(
\sigma\left(  i\right)  -i\right) \\
&  =\underbrace{\sum_{\substack{i\in\left[  n\right]  ;\\\sigma\left(
i\right)  <i}}\left(  i-\sigma\left(  i\right)  \right)  }_{\substack{=\sum
_{\substack{i\in\left[  n\right]  ;\\\sigma\left(  i\right)  >i}}\left(
\sigma\left(  i\right)  -i\right)  \\\text{(by
(\ref{sol.exe.perm.footrule.a.1}))}}}+\underbrace{\sum_{\substack{i\in\left[
n\right]  ;\\\sigma\left(  i\right)  =i}}0}_{=0}+\sum_{\substack{i\in\left[
n\right]  ;\\\sigma\left(  i\right)  >i}}\left(  \sigma\left(  i\right)
-i\right) \\
&  =\sum_{\substack{i\in\left[  n\right]  ;\\\sigma\left(  i\right)
>i}}\left(  \sigma\left(  i\right)  -i\right)  +\sum_{\substack{i\in\left[
n\right]  ;\\\sigma\left(  i\right)  >i}}\left(  \sigma\left(  i\right)
-i\right)  =2\sum_{\substack{i\in\left[  n\right]  ;\\\sigma\left(  i\right)
>i}}\left(  \sigma\left(  i\right)  -i\right)  .
\end{align*}
Now, the definition of $h\left(  \sigma\right)  $ yields%
\[
h\left(  \sigma\right)  =\sum_{i\in\left[  n\right]  }\left\vert \sigma\left(
i\right)  -i\right\vert =2\underbrace{\sum_{\substack{i\in\left[  n\right]
;\\\sigma\left(  i\right)  >i}}\left(  \sigma\left(  i\right)  -i\right)
}_{\substack{=\sum_{\substack{i\in\left[  n\right]  ;\\\sigma\left(  i\right)
<i}}\left(  i-\sigma\left(  i\right)  \right)  \\\text{(by
(\ref{sol.exe.perm.footrule.a.1}))}}}=2\sum_{\substack{i\in\left[  n\right]
;\\\sigma\left(  i\right)  <i}}\left(  i-\sigma\left(  i\right)  \right)  .
\]
This solves Exercise \ref{exe.perm.footrule} \textbf{(a)}.

\textbf{(b)} Let $\tau\in S_{n}$.

\begin{vershort}
We have proven that $\sigma$ is a bijection $\left[  n\right]  \rightarrow
\left[  n\right]  $. Similarly, $\tau$ is a bijection $\left[  n\right]
\rightarrow\left[  n\right]  $.
\end{vershort}

\begin{verlong}
We have proven that $\sigma$ is a bijection $\left[  n\right]  \rightarrow
\left[  n\right]  $. The same argument (applied to $\tau$ instead of $\sigma$)
shows that $\tau$ is a bijection $\left[  n\right]  \rightarrow\left[
n\right]  $.
\end{verlong}

The definition of $h\left(  \tau\right)  $ yields%
\begin{equation}
h\left(  \tau\right)  =\sum_{i\in\left[  n\right]  }\left\vert \tau\left(
i\right)  -i\right\vert . \label{sol.perm.footrule.b.1}%
\end{equation}
The definition of $h\left(  \sigma\circ\tau\right)  $ yields%
\begin{equation}
h\left(  \sigma\circ\tau\right)  =\sum_{i\in\left[  n\right]  }\left\vert
\underbrace{\left(  \sigma\circ\tau\right)  \left(  i\right)  }_{=\sigma
\left(  \tau\left(  i\right)  \right)  }-i\right\vert =\sum_{i\in\left[
n\right]  }\left\vert \sigma\left(  \tau\left(  i\right)  \right)
-i\right\vert . \label{sol.perm.footrule.b.2}%
\end{equation}
But the definition of $h\left(  \sigma\right)  $ yields%
\[
h\left(  \sigma\right)  =\sum_{i\in\left[  n\right]  }\left\vert \sigma\left(
i\right)  -i\right\vert =\sum_{i\in\left[  n\right]  }\left\vert \sigma\left(
\tau\left(  i\right)  \right)  -\tau\left(  i\right)  \right\vert
\]
(here, we have substituted $\tau\left(  i\right)  $ for $i$ in the sum, since
$\tau$ is a bijection $\left[  n\right]  \rightarrow\left[  n\right]  $).
Adding (\ref{sol.perm.footrule.b.1}) to this equality, we obtain%
\begin{align*}
h\left(  \sigma\right)  +h\left(  \tau\right)   &  =\sum_{i\in\left[
n\right]  }\left\vert \sigma\left(  \tau\left(  i\right)  \right)
-\tau\left(  i\right)  \right\vert +\sum_{i\in\left[  n\right]  }\left\vert
\tau\left(  i\right)  -i\right\vert \\
&  =\sum_{i\in\left[  n\right]  }\underbrace{\left(  \left\vert \sigma\left(
\tau\left(  i\right)  \right)  -\tau\left(  i\right)  \right\vert +\left\vert
\tau\left(  i\right)  -i\right\vert \right)  }_{\substack{\geq\left\vert
\left(  \sigma\left(  \tau\left(  i\right)  \right)  -\tau\left(  i\right)
\right)  +\left(  \tau\left(  i\right)  -i\right)  \right\vert \\\text{(by
Theorem \ref{thm.ineq.triangle-R}}\\\text{(applied to }x=\sigma\left(
\tau\left(  i\right)  \right)  -\tau\left(  i\right)  \text{ and }%
y=\tau\left(  i\right)  -i\text{))}}}\\
&  \geq\sum_{i\in\left[  n\right]  }\left\vert \underbrace{\left(
\sigma\left(  \tau\left(  i\right)  \right)  -\tau\left(  i\right)  \right)
+\left(  \tau\left(  i\right)  -i\right)  }_{=\sigma\left(  \tau\left(
i\right)  \right)  -i}\right\vert =\sum_{i\in\left[  n\right]  }\left\vert
\sigma\left(  \tau\left(  i\right)  \right)  -i\right\vert =h\left(
\sigma\circ\tau\right)
\end{align*}
(by (\ref{sol.perm.footrule.b.2})). In other words, $h\left(  \sigma\circ
\tau\right)  \leq h\left(  \sigma\right)  +h\left(  \tau\right)  $. This
solves Exercise \ref{exe.perm.footrule} \textbf{(b)}.

\textbf{(d)} Exercise \ref{exe.perm.footrule} \textbf{(a)} yields $h\left(
\sigma\right)  =2\sum_{\substack{i\in\left[  n\right]  ;\\\sigma\left(
i\right)  >i}}\left(  \sigma\left(  i\right)  -i\right)  $. Dividing this
equality by $2$, we obtain%
\begin{align}
h\left(  \sigma\right)  /2  &  =\sum_{\substack{i\in\left[  n\right]
;\\\sigma\left(  i\right)  >i}}\left(  \underbrace{\sigma\left(  i\right)
}_{\substack{\leq i+\ell_{i}\left(  \sigma\right)  \\\text{(by Lemma
\ref{lem.perm.lexico1.lis} \textbf{(c)})}}}-i\right)
\label{sol.perm.footrule.d.l.0a}\\
&  \leq\sum_{\substack{i\in\left[  n\right]  ;\\\sigma\left(  i\right)
>i}}\underbrace{\left(  i+\ell_{i}\left(  \sigma\right)  -i\right)  }%
_{=\ell_{i}\left(  \sigma\right)  }=\sum_{\substack{i\in\left[  n\right]
;\\\sigma\left(  i\right)  >i}}\ell_{i}\left(  \sigma\right)  .
\label{sol.perm.footrule.d.l.1}%
\end{align}
Also, Exercise \ref{exe.perm.footrule} \textbf{(a)} yields $h\left(
\sigma\right)  =2\sum_{\substack{i\in\left[  n\right]  ;\\\sigma\left(
i\right)  <i}}\left(  i-\sigma\left(  i\right)  \right)  $. Dividing this
equality by $2$, we obtain%
\begin{equation}
h\left(  \sigma\right)  /2=\sum_{\substack{i\in\left[  n\right]
;\\\sigma\left(  i\right)  <i}}\left(  i-\sigma\left(  i\right)  \right)  .
\label{sol.perm.footrule.d.l.0b}%
\end{equation}

On the other hand, each $i\in\left[  n\right]  $ satisfies
\begin{equation}
\ell_{i}\left(  \sigma\right)  \geq0. \label{sol.perm.footrule.d.l.geq0}%
\end{equation}

\begin{vershort}
\noindent(This follows immediately from the definition of $\ell_{i}\left(
\sigma\right)  $.)
\end{vershort}

\begin{verlong}
[\textit{Proof of (\ref{sol.perm.footrule.d.l.geq0}):} Let $i\in\left[
n\right]  $. Then, $\ell_{i}\left(  \sigma\right)  $ is the number of all
$j\in\left\{  i+1,i+2,\ldots,n\right\}  $ such that $\sigma\left(  i\right)
>\sigma\left(  j\right)  $ (by the definition of $\ell_{i}\left(
\sigma\right)  $). Thus, $\ell_{i}\left(  \sigma\right)  $ is a nonnegative
integer. Hence, $\ell_{i}\left(  \sigma\right)  \geq0$. This proves
(\ref{sol.perm.footrule.d.l.geq0}).]
\end{verlong}

Proposition \ref{prop.perm.lehmer.l} yields%
\begin{align}
\ell\left(  \sigma\right)   &  =\ell_{1}\left(  \sigma\right)  +\ell
_{2}\left(  \sigma\right)  +\cdots+\ell_{n}\left(  \sigma\right)
=\underbrace{\sum_{i=1}^{n}}_{\substack{=\sum_{i\in\left\{  1,2,\ldots
,n\right\}  }=\sum_{i\in\left[  n\right]  }\\\text{(since }\left\{
1,2,\ldots,n\right\}  =\left[  n\right]  \text{)}}}\ell_{i}\left(
\sigma\right) \nonumber\\
&  =\sum_{i\in\left[  n\right]  }\ell_{i}\left(  \sigma\right)
\label{sol.perm.footrule.d.lsig=1}\\
&  =\sum_{\substack{i\in\left[  n\right]  ;\\\sigma\left(  i\right)  \leq
i}}\underbrace{\ell_{i}\left(  \sigma\right)  }_{\substack{\geq0\\\text{(by
(\ref{sol.perm.footrule.d.l.geq0}))}}}+\underbrace{\sum_{\substack{i\in\left[
n\right]  ;\\\sigma\left(  i\right)  >i}}\ell_{i}\left(  \sigma\right)
}_{\substack{\geq h\left(  \sigma\right)  /2\\\text{(by
(\ref{sol.perm.footrule.d.l.1}))}}}\nonumber\\
&  \ \ \ \ \ \ \ \ \ \ \left(
\begin{array}
[c]{c}%
\text{since each }i\in\left[  n\right]  \text{ satisfies either }\sigma\left(
i\right)  \leq i\text{ or }\sigma\left(  i\right)  >i\text{,}\\
\text{but not both at the same time}%
\end{array}
\right) \nonumber\\
&  \geq\underbrace{\sum_{\substack{i\in\left[  n\right]  ;\\\sigma\left(
i\right)  \leq i}}0}_{=0}+h\left(  \sigma\right)  /2\geq h\left(
\sigma\right)  /2.\nonumber
\end{align}
In other words, $h\left(  \sigma\right)  /2\leq\ell\left(  \sigma\right)  $.

It now remains to prove that $\ell\left(  \sigma\right)  \leq h\left(
\sigma\right)  $. Our proof of this inequality shall follow the idea given in
the Hint, but not literally: Instead of counting inversions, we will be
summing truth values (which boils down to the same thing, but looks slicker on paper).

The equality (\ref{sol.perm.footrule.d.lsig=1}) becomes%
\begin{align}
\ell\left(  \sigma\right)   &  =\sum_{i\in\left[  n\right]  }\underbrace{\ell
_{i}\left(  \sigma\right)  }_{\substack{=\sum_{j\in\left[  n\right]  }\left[
i<j\right]  \left[  \sigma\left(  i\right)  >\sigma\left(  j\right)  \right]
\\\text{(by Lemma \ref{lem.sol.perm.lisitau.1} \textbf{(b)})}}}\nonumber\\
&  =\sum_{i\in\left[  n\right]  }\sum_{j\in\left[  n\right]  }\left[
i<j\right]  \left[  \sigma\left(  i\right)  >\sigma\left(  j\right)  \right]
. \label{sol.perm.footrule.d.r.1}%
\end{align}

Next, we claim:

\begin{statement}
\textit{Claim 1:} Let $i\in\left[  n\right]  $ and $j\in\left[  n\right]  $.
Then,%
\[
\left[  i<j\right]  \left[  \sigma\left(  i\right)  >\sigma\left(  j\right)
\right]  \leq\left[  \sigma\left(  i\right)  >\sigma\left(  j\right)  \right]
\left[  \sigma\left(  i\right)  <j\right]  +\left[  j>i\right]  \left[
j\leq\sigma\left(  i\right)  \right]  .
\]

\end{statement}

\begin{vershort}
[\textit{Proof of Claim 1:} The logical statements $\left(  i<j\right)  $ and
$\left(  j>i\right)  $ are equivalent. Hence, Exercise \ref{exe.iverson-prop}
\textbf{(a)} yields $\left[  i<j\right]  =\left[  j>i\right]  $.

We have $\sigma\left(  i\right)  \in\left[  n\right]  $ (since $\sigma$ is a
map $\left[  n\right]  \rightarrow\left[  n\right]  $). Now,%
\begin{align*}
&  \left[  i<j\right]  \left[  \sigma\left(  i\right)  >\sigma\left(
j\right)  \right] \\
&  =\left[  i<j\right]  \left[  \sigma\left(  i\right)  >\sigma\left(
j\right)  \right]  \cdot\underbrace{1}_{\substack{=\left[  \sigma\left(
i\right)  <j\right]  +\left[  j\leq\sigma\left(  i\right)  \right]
\\\text{(by Lemma \ref{lem.sol.perm.footrule.1} \textbf{(a)}}\\\text{(applied
to }u=\sigma\left(  i\right)  \text{ and }v=j\text{))}}}\\
&  =\left[  i<j\right]  \left[  \sigma\left(  i\right)  >\sigma\left(
j\right)  \right]  \left(  \left[  \sigma\left(  i\right)  <j\right]  +\left[
j\leq\sigma\left(  i\right)  \right]  \right) \\
&  =\underbrace{\left[  i<j\right]  }_{\substack{\leq1\\\text{(since }\left[
\mathcal{A}\right]  \leq1\text{ for}\\\text{any statement }\mathcal{A}%
\text{)}}}\left[  \sigma\left(  i\right)  >\sigma\left(  j\right)  \right]
\left[  \sigma\left(  i\right)  <j\right]  +\underbrace{\left[  i<j\right]
}_{=\left[  j>i\right]  }\underbrace{\left[  \sigma\left(  i\right)
>\sigma\left(  j\right)  \right]  }_{\substack{\leq1\\\text{(since }\left[
\mathcal{A}\right]  \leq1\text{ for}\\\text{any statement }\mathcal{A}%
\text{)}}}\left[  j\leq\sigma\left(  i\right)  \right] \\
&  \leq\left[  \sigma\left(  i\right)  >\sigma\left(  j\right)  \right]
\left[  \sigma\left(  i\right)  <j\right]  +\left[  j>i\right]  \left[
j\leq\sigma\left(  i\right)  \right]  .
\end{align*}
(Here, we have made tacit use of the fact that all terms in our computation
are nonnegative\footnote{because $\left[  \mathcal{A}\right]  \geq0$ for any
statement $\mathcal{A}$}. Indeed, this fact allowed us to multiply
inequalities.) Thus, Claim 1 is proven.]
\end{vershort}

\begin{verlong}
[\textit{Proof of Claim 1:} The logical statements $\left(  i<j\right)  $ and
$\left(  j>i\right)  $ are equivalent. Hence, Exercise \ref{exe.iverson-prop}
\textbf{(a)} (applied to $\mathcal{A}=\left(  i<j\right)  $ and $\mathcal{B}%
=\left(  j>i\right)  $) yields $\left[  i<j\right]  =\left[  j>i\right]  $.

Every logical statement $\mathcal{A}$ satisfies $\left[  \mathcal{A}\right]
\in\left\{  0,1\right\}  $ (by the definition of $\left[  \mathcal{A}\right]
$) and thus $\left[  \mathcal{A}\right]  \leq1$. Hence, $\left[  i<j\right]
\leq1$ and $\left[  \sigma\left(  i\right)  >\sigma\left(  j\right)  \right]
\leq1$.

Every logical statement $\mathcal{A}$ satisfies $\left[  \mathcal{A}\right]
\in\left\{  0,1\right\}  $ (by the definition of $\left[  \mathcal{A}\right]
$) and thus $\left[  \mathcal{A}\right]  \geq0$. Hence, $\left[  \sigma\left(
i\right)  >\sigma\left(  j\right)  \right]  \geq0$ and $\left[  \sigma\left(
i\right)  <j\right]  \geq0$ and $\left[  i<j\right]  \geq0$ and $\left[
j\leq\sigma\left(  i\right)  \right]  \geq0$.

The number $\left[  \sigma\left(  i\right)  >\sigma\left(  j\right)  \right]
\left[  \sigma\left(  i\right)  <j\right]  $ is nonnegative (since $\left[
\sigma\left(  i\right)  >\sigma\left(  j\right)  \right]  \geq0$ and $\left[
\sigma\left(  i\right)  <j\right]  \geq0$). Thus, we can multiply the
inequality $\left[  i<j\right]  \leq1$ by this number. We thus obtain%
\begin{align}
\left[  i<j\right]  \left[  \sigma\left(  i\right)  >\sigma\left(  j\right)
\right]  \left[  \sigma\left(  i\right)  <j\right]   &  \leq1\left[
\sigma\left(  i\right)  >\sigma\left(  j\right)  \right]  \left[
\sigma\left(  i\right)  <j\right] \nonumber\\
&  =\left[  \sigma\left(  i\right)  >\sigma\left(  j\right)  \right]  \left[
\sigma\left(  i\right)  <j\right]  . \label{sol.perm.footrule.d.c1.pf.1}%
\end{align}

The number $\left[  i<j\right]  \left[  j\leq\sigma\left(  i\right)  \right]
$ is nonnegative (since $\left[  i<j\right]  \geq0$ and $\left[  j\leq
\sigma\left(  i\right)  \right]  \geq0$). Thus, we can multiply the inequality
$\left[  \sigma\left(  i\right)  >\sigma\left(  j\right)  \right]  \leq1$ by
this number. We thus obtain%
\begin{align}
\left[  \sigma\left(  i\right)  >\sigma\left(  j\right)  \right]  \left[
i<j\right]  \left[  j\leq\sigma\left(  i\right)  \right]   &  \leq1\left[
i<j\right]  \left[  j\leq\sigma\left(  i\right)  \right]  =\underbrace{\left[
i<j\right]  }_{=\left[  j>i\right]  }\left[  j\leq\sigma\left(  i\right)
\right] \nonumber\\
&  =\left[  j>i\right]  \left[  j\leq\sigma\left(  i\right)  \right]  .
\label{sol.perm.footrule.d.c1.pf.2}%
\end{align}

We have $\sigma\left(  i\right)  \in\left[  n\right]  $ (since $\sigma$ is a
map $\left[  n\right]  \rightarrow\left[  n\right]  $). Now,%
\begin{align*}
&  \left[  i<j\right]  \left[  \sigma\left(  i\right)  >\sigma\left(
j\right)  \right] \\
&  =\left[  i<j\right]  \left[  \sigma\left(  i\right)  >\sigma\left(
j\right)  \right]  \cdot\underbrace{1}_{\substack{=\left[  \sigma\left(
i\right)  <j\right]  +\left[  j\leq\sigma\left(  i\right)  \right]
\\\text{(by Lemma \ref{lem.sol.perm.footrule.1} \textbf{(a)}}\\\text{(applied
to }u=\sigma\left(  i\right)  \text{ and }v=j\text{))}}}\\
&  =\left[  i<j\right]  \left[  \sigma\left(  i\right)  >\sigma\left(
j\right)  \right]  \left(  \left[  \sigma\left(  i\right)  <j\right]  +\left[
j\leq\sigma\left(  i\right)  \right]  \right) \\
&  =\underbrace{\left[  i<j\right]  \left[  \sigma\left(  i\right)
>\sigma\left(  j\right)  \right]  \left[  \sigma\left(  i\right)  <j\right]
}_{\substack{\leq\left[  \sigma\left(  i\right)  >\sigma\left(  j\right)
\right]  \left[  \sigma\left(  i\right)  <j\right]  \\\text{(by
(\ref{sol.perm.footrule.d.c1.pf.1}))}}}+\underbrace{\left[  i<j\right]
\left[  \sigma\left(  i\right)  >\sigma\left(  j\right)  \right]  \left[
j\leq\sigma\left(  i\right)  \right]  }_{\substack{=\left[  \sigma\left(
i\right)  >\sigma\left(  j\right)  \right]  \left[  i<j\right]  \left[
j\leq\sigma\left(  i\right)  \right]  \\\leq\left[  j>i\right]  \left[
j\leq\sigma\left(  i\right)  \right]  \\\text{(by
(\ref{sol.perm.footrule.d.c1.pf.2}))}}}\\
&  \leq\left[  \sigma\left(  i\right)  >\sigma\left(  j\right)  \right]
\left[  \sigma\left(  i\right)  <j\right]  +\left[  j>i\right]  \left[
j\leq\sigma\left(  i\right)  \right]  .
\end{align*}
This proves Claim 1.]
\end{verlong}

Now, (\ref{sol.perm.footrule.d.r.1}) becomes%
\begin{align}
\ell\left(  \sigma\right)   &  =\sum_{i\in\left[  n\right]  }\sum_{j\in\left[
n\right]  }\underbrace{\left[  i<j\right]  \left[  \sigma\left(  i\right)
>\sigma\left(  j\right)  \right]  }_{\substack{\leq\left[  \sigma\left(
i\right)  >\sigma\left(  j\right)  \right]  \left[  \sigma\left(  i\right)
<j\right]  +\left[  j>i\right]  \left[  j\leq\sigma\left(  i\right)  \right]
\\\text{(by Claim 1)}}}\nonumber\\
&  \leq\sum_{i\in\left[  n\right]  }\sum_{j\in\left[  n\right]  }\left(
\left[  \sigma\left(  i\right)  >\sigma\left(  j\right)  \right]  \left[
\sigma\left(  i\right)  <j\right]  +\left[  j>i\right]  \left[  j\leq
\sigma\left(  i\right)  \right]  \right) \nonumber\\
&  =\underbrace{\sum_{i\in\left[  n\right]  }\sum_{j\in\left[  n\right]  }%
}_{=\sum_{j\in\left[  n\right]  }\sum_{i\in\left[  n\right]  }}\left[
\sigma\left(  i\right)  >\sigma\left(  j\right)  \right]  \left[
\sigma\left(  i\right)  <j\right]  +\sum_{i\in\left[  n\right]  }\sum
_{j\in\left[  n\right]  }\left[  j>i\right]  \left[  j\leq\sigma\left(
i\right)  \right] \nonumber\\
&  =\sum_{j\in\left[  n\right]  }\underbrace{\sum_{i\in\left[  n\right]
}\left[  \sigma\left(  i\right)  >\sigma\left(  j\right)  \right]  \left[
\sigma\left(  i\right)  <j\right]  }_{\substack{=\sum_{k\in\left[  n\right]
}\left[  k>\sigma\left(  j\right)  \right]  \left[  k<j\right]  \\\text{(here,
we have substituted }k\text{ for }\sigma\left(  i\right)  \text{
in}\\\text{the sum, since the map }\sigma:\left[  n\right]  \rightarrow\left[
n\right]  \text{ is a bijection)}}}+\sum_{i\in\left[  n\right]  }%
\underbrace{\sum_{j\in\left[  n\right]  }\left[  j>i\right]  \left[
j\leq\sigma\left(  i\right)  \right]  }_{\substack{=\sum_{k\in\left[
n\right]  }\left[  k>i\right]  \left[  k\leq\sigma\left(  i\right)  \right]
\\\text{(here, we have renamed the summation}\\\text{index }j\text{ as
}k\text{)}}}\nonumber\\
&  =\underbrace{\sum_{j\in\left[  n\right]  }\sum_{k\in\left[  n\right]
}\left[  k>\sigma\left(  j\right)  \right]  \left[  k<j\right]  }%
_{\substack{=\sum_{i\in\left[  n\right]  }\sum_{k\in\left[  n\right]  }\left[
k>\sigma\left(  i\right)  \right]  \left[  k<i\right]  \\\text{(here, we have
renamed the summation index }j\\\text{as }i\text{ in the outer sum)}}%
}+\sum_{i\in\left[  n\right]  }\sum_{k\in\left[  n\right]  }\left[
k>i\right]  \left[  k\leq\sigma\left(  i\right)  \right] \nonumber\\
&  =\sum_{i\in\left[  n\right]  }\sum_{k\in\left[  n\right]  }\left[
k>\sigma\left(  i\right)  \right]  \left[  k<i\right]  +\sum_{i\in\left[
n\right]  }\sum_{k\in\left[  n\right]  }\left[  k>i\right]  \left[
k\leq\sigma\left(  i\right)  \right]  . \label{sol.perm.footrule.d.r.3}%
\end{align}

Now, fix $i\in\left[  n\right]  $. Then, $\sigma\left(  i\right)  \in\left[
n\right]  $ (since $\sigma$ is a map $\left[  n\right]  \rightarrow\left[
n\right]  $). Hence, Lemma \ref{lem.sol.perm.footrule.1} \textbf{(b)} (applied
to $u=\sigma\left(  i\right)  $ and $v=i$) yields%
\begin{equation}
\sum_{k\in\left[  n\right]  }\left[  k>\sigma\left(  i\right)  \right]
\left[  k<i\right]  =\left[  \sigma\left(  i\right)  <i\right]  \left(
i-\sigma\left(  i\right)  -1\right)  . \label{sol.perm.footrule.d.r.4a}%
\end{equation}
Moreover, Lemma \ref{lem.sol.perm.footrule.1} \textbf{(c)} (applied to $u=i$
and $v=\sigma\left(  i\right)  $) yields%
\begin{equation}
\sum_{k\in\left[  n\right]  }\left[  k>i\right]  \left[  k\leq\sigma\left(
i\right)  \right]  =\left[  i<\sigma\left(  i\right)  \right]  \left(
\sigma\left(  i\right)  -i\right)  . \label{sol.perm.footrule.d.r.4b}%
\end{equation}

Now, forget that we have fixed $i$. We thus have proven
(\ref{sol.perm.footrule.d.r.4a}) and (\ref{sol.perm.footrule.d.r.4b}) for each
$i\in\left[  n\right]  $.

Now,
\begin{align}
&  \sum_{i\in\left[  n\right]  }\underbrace{\sum_{k\in\left[  n\right]
}\left[  k>\sigma\left(  i\right)  \right]  \left[  k<i\right]  }%
_{\substack{=\left[  \sigma\left(  i\right)  <i\right]  \left(  i-\sigma
\left(  i\right)  -1\right)  \\\text{(by (\ref{sol.perm.footrule.d.r.4a}))}%
}}\nonumber\\
&  =\sum_{i\in\left[  n\right]  }\left[  \sigma\left(  i\right)  <i\right]
\left(  i-\sigma\left(  i\right)  -1\right) \nonumber\\
&  =\sum_{\substack{i\in\left[  n\right]  ;\\\sigma\left(  i\right)
<i}}\underbrace{\left[  \sigma\left(  i\right)  <i\right]  }%
_{\substack{=1\\\text{(since }\sigma\left(  i\right)  <i\text{)}}}\left(
i-\sigma\left(  i\right)  -1\right)  +\sum_{\substack{i\in\left[  n\right]
;\\\text{not }\sigma\left(  i\right)  <i}}\underbrace{\left[  \sigma\left(
i\right)  <i\right]  }_{\substack{=0\\\text{(since we don't have }%
\sigma\left(  i\right)  <i\text{)}}}\left(  i-\sigma\left(  i\right)
-1\right) \nonumber\\
&  \ \ \ \ \ \ \ \ \ \ \left(
\begin{array}
[c]{c}%
\text{since each }i\in\left[  n\right]  \text{ satisfies either }\sigma\left(
i\right)  <i\text{ or }\left(  \text{not }\sigma\left(  i\right)  <i\right)
\text{,}\\
\text{but not both at the same time}%
\end{array}
\right) \nonumber\\
&  =\sum_{\substack{i\in\left[  n\right]  ;\\\sigma\left(  i\right)
<i}}\left(  i-\sigma\left(  i\right)  -1\right)  +\underbrace{\sum
_{\substack{i\in\left[  n\right]  ;\\\text{not }\sigma\left(  i\right)
<i}}0\left(  i-\sigma\left(  i\right)  -1\right)  }_{=0}=\sum_{\substack{i\in
\left[  n\right]  ;\\\sigma\left(  i\right)  <i}}\underbrace{\left(
i-\sigma\left(  i\right)  -1\right)  }_{\leq i-\sigma\left(  i\right)
}\nonumber\\
&  \leq\sum_{\substack{i\in\left[  n\right]  ;\\\sigma\left(  i\right)
<i}}\left(  i-\sigma\left(  i\right)  \right)  =h\left(  \sigma\right)
/2\ \ \ \ \ \ \ \ \ \ \left(  \text{by (\ref{sol.perm.footrule.d.l.0b}%
)}\right)  . \label{sol.perm.footrule.d.r.5a}%
\end{align}

Furthermore,%
\begin{align}
&  \sum_{i\in\left[  n\right]  }\underbrace{\sum_{k\in\left[  n\right]
}\left[  k>i\right]  \left[  k\leq\sigma\left(  i\right)  \right]
}_{\substack{=\left[  i<\sigma\left(  i\right)  \right]  \left(  \sigma\left(
i\right)  -i\right)  \\\text{(by (\ref{sol.perm.footrule.d.r.4b}))}%
}}\nonumber\\
&  =\sum_{i\in\left[  n\right]  }\left[  i<\sigma\left(  i\right)  \right]
\left(  \sigma\left(  i\right)  -i\right) \nonumber\\
&  =\sum_{\substack{i\in\left[  n\right]  ;\\i<\sigma\left(  i\right)
}}\underbrace{\left[  i<\sigma\left(  i\right)  \right]  }%
_{\substack{=1\\\text{(since }i<\sigma\left(  i\right)  \text{)}}}\left(
\sigma\left(  i\right)  -i\right)  +\sum_{\substack{i\in\left[  n\right]
;\\\text{not }i<\sigma\left(  i\right)  }}\underbrace{\left[  i<\sigma\left(
i\right)  \right]  }_{\substack{=0\\\text{(since we don't have }%
i<\sigma\left(  i\right)  \text{)}}}\left(  \sigma\left(  i\right)  -i\right)
\nonumber\\
&  \ \ \ \ \ \ \ \ \ \ \left(
\begin{array}
[c]{c}%
\text{since each }i\in\left[  n\right]  \text{ satisfies either }%
i<\sigma\left(  i\right)  \text{ or }\left(  \text{not }i<\sigma\left(
i\right)  \right)  \text{,}\\
\text{but not both at the same time}%
\end{array}
\right) \nonumber\\
&  =\sum_{\substack{i\in\left[  n\right]  ;\\i<\sigma\left(  i\right)
}}\left(  \sigma\left(  i\right)  -i\right)  +\underbrace{\sum_{\substack{i\in
\left[  n\right]  ;\\\text{not }i<\sigma\left(  i\right)  }}0\left(
\sigma\left(  i\right)  -i\right)  }_{=0}=\underbrace{\sum_{\substack{i\in
\left[  n\right]  ;\\i<\sigma\left(  i\right)  }}}_{\substack{=\sum
_{\substack{i\in\left[  n\right]  ;\\\sigma\left(  i\right)  >i}%
}\\\text{(since the condition }\left(  i<\sigma\left(  i\right)  \right)
\\\text{is equivalent to }\left(  \sigma\left(  i\right)  >i\right)  \text{)}%
}}\left(  \sigma\left(  i\right)  -i\right) \nonumber\\
&  =\sum_{\substack{i\in\left[  n\right]  ;\\\sigma\left(  i\right)
>i}}\left(  \sigma\left(  i\right)  -i\right)  =h\left(  \sigma\right)
/2\ \ \ \ \ \ \ \ \ \ \left(  \text{by (\ref{sol.perm.footrule.d.l.0a}%
)}\right)  . \label{sol.perm.footrule.d.r.5b}%
\end{align}

Now, (\ref{sol.perm.footrule.d.r.3}) becomes%
\begin{align*}
\ell\left(  \sigma\right)   &  \leq\underbrace{\sum_{i\in\left[  n\right]
}\sum_{k\in\left[  n\right]  }\left[  k>\sigma\left(  i\right)  \right]
\left[  k<i\right]  }_{\substack{\leq h\left(  \sigma\right)  /2\\\text{(by
(\ref{sol.perm.footrule.d.r.5a}))}}}+\underbrace{\sum_{i\in\left[  n\right]
}\sum_{k\in\left[  n\right]  }\left[  k>i\right]  \left[  k\leq\sigma\left(
i\right)  \right]  }_{\substack{=h\left(  \sigma\right)  /2\\\text{(by
(\ref{sol.perm.footrule.d.r.5b}))}}}\\
&  \leq h\left(  \sigma\right)  /2+h\left(  \sigma\right)  /2=h\left(
\sigma\right)  .
\end{align*}
Thus, $\ell\left(  \sigma\right)  \leq h\left(  \sigma\right)  $ is proven.

We have now proven both $h\left(  \sigma\right)  /2\leq\ell\left(
\sigma\right)  $ and $\ell\left(  \sigma\right)  \leq h\left(  \sigma\right)
$. Thus, $h\left(  \sigma\right)  /2\leq\ell\left(  \sigma\right)  \leq
h\left(  \sigma\right)  $ follows. This solves Exercise
\ref{exe.perm.footrule} \textbf{(d)}.

\textbf{(c)} Let $k\in\left\{  1,2,\ldots,n-1\right\}  $. We have $\ell\left(
s_{k}\right)  =1$. (This was shown during our proof of Proposition
\ref{prop.perm.signs.basics} \textbf{(b)}.) But Exercise
\ref{exe.perm.footrule} \textbf{(d)} (applied to $s_{k}$ instead of $\sigma$)
yields $h\left(  s_{k}\right)  /2\leq\ell\left(  s_{k}\right)  \leq h\left(
s_{k}\right)  $. From $h\left(  s_{k}\right)  /2\leq\ell\left(  s_{k}\right)
$, we obtain $h\left(  s_{k}\right)  \leq2\underbrace{\ell\left(
s_{k}\right)  }_{=1}=2$. (It is not hard to check that $h\left(  s_{k}\right)
=2$, but we will not need this fact.)

Now, Exercise \ref{exe.perm.footrule} \textbf{(b)} (applied to $s_{k}$ and
$\sigma$ instead of $\sigma$ and $\tau$) yields $h\left(  s_{k}\circ
\sigma\right)  \leq\underbrace{h\left(  s_{k}\right)  }_{\leq2}+h\left(
\sigma\right)  \leq2+h\left(  \sigma\right)  =h\left(  \sigma\right)  +2$.
This solves Exercise \ref{exe.perm.footrule} \textbf{(c)}.
\end{proof}

\subsection{Solution to Exercise \ref{exe.perm.Inv.sub}}

We shall now prepare for the solution of Exercise \ref{exe.perm.Inv.sub}.
First, we introduce a notation:

\begin{definition}
\label{def.sol.perm.Inv.sub.aXb}If $X$, $X^{\prime}$, $Y$ and $Y^{\prime}$ are
four sets and if $\alpha:X\rightarrow X^{\prime}$ and $\beta:Y\rightarrow
Y^{\prime}$ are two maps, then $\alpha\times\beta$ will denote the map%
\begin{align*}
X\times Y  &  \rightarrow X^{\prime}\times Y^{\prime},\\
\left(  x,y\right)   &  \mapsto\left(  \alpha\left(  x\right)  ,\beta\left(
y\right)  \right)  .
\end{align*}

\end{definition}

\begin{lemma}
\label{lem.sol.perm.Inv.sub.bijbij}Let $n\in\mathbb{N}$. Let $\left[
n\right]  $ denote the set $\left\{  1,2,\ldots,n\right\}  $. Let $\alpha\in
S_{n}$ and $\beta\in S_{n}$. Then, the map $\alpha\times\beta:\left[
n\right]  \times\left[  n\right]  \rightarrow\left[  n\right]  \times\left[
n\right]  $ is invertible, and its inverse is $\alpha^{-1}\times\beta^{-1}$.
\end{lemma}

\begin{proof}
[Proof of Lemma \ref{lem.sol.perm.Inv.sub.bijbij}.]Recall that $S_{n}$ is the
set of all permutations of the set $\left\{  1,2,\ldots,n\right\}  $. In other
words, $S_{n}$ is the set of all permutations of the set $\left[  n\right]  $
(since $\left\{  1,2,\ldots,n\right\}  =\left[  n\right]  $).

We know that $\alpha\in S_{n}$. In other words, $\alpha$ is a permutation of
the set $\left[  n\right]  $ (since $S_{n}$ is the set of all permutations of
the set $\left[  n\right]  $). In other words, $\alpha$ is a bijective map
$\left[  n\right]  \rightarrow\left[  n\right]  $. Similarly, $\beta$ is a
bijective map $\left[  n\right]  \rightarrow\left[  n\right]  $. Hence,
Definition \ref{def.sol.perm.Inv.sub.aXb} defines a map $\alpha\times
\beta:\left[  n\right]  \times\left[  n\right]  \rightarrow\left[  n\right]
\times\left[  n\right]  $.

The map $\alpha$ is bijective, and thus invertible. Hence, its inverse map
$\alpha^{-1}:\left[  n\right]  \rightarrow\left[  n\right]  $ is well-defined.
Similarly, $\beta^{-1}:\left[  n\right]  \rightarrow\left[  n\right]  $ is
well-defined. Thus, Definition \ref{def.sol.perm.Inv.sub.aXb} defines a map
$\alpha^{-1}\times\beta^{-1}:\left[  n\right]  \times\left[  n\right]
\rightarrow\left[  n\right]  \times\left[  n\right]  $.

\begin{vershort}
Now, every $\left(  i,j\right)  \in\left[  n\right]  \times\left[  n\right]  $
satisfies
\begin{align*}
&  \left(  \left(  \alpha^{-1}\times\beta^{-1}\right)  \circ\left(
\alpha\times\beta\right)  \right)  \left(  i,j\right) \\
&  =\left(  \alpha^{-1}\times\beta^{-1}\right)  \left(  \underbrace{\left(
\alpha\times\beta\right)  \left(  i,j\right)  }_{\substack{=\left(
\alpha\left(  i\right)  ,\beta\left(  j\right)  \right)  \\\text{(by the
definition of }\alpha\times\beta\text{)}}}\right) \\
&  =\left(  \alpha^{-1}\times\beta^{-1}\right)  \left(  \alpha\left(
i\right)  ,\beta\left(  j\right)  \right)  =\left(  \underbrace{\alpha
^{-1}\left(  \alpha\left(  i\right)  \right)  }_{=i},\underbrace{\beta
^{-1}\left(  \beta\left(  j\right)  \right)  }_{=j}\right) \\
&  \ \ \ \ \ \ \ \ \ \ \left(  \text{by the definition of }\alpha^{-1}%
\times\beta^{-1}\right) \\
&  =\left(  i,j\right)  =\operatorname*{id}\left(  i,j\right)  .
\end{align*}
Thus, $\left(  \alpha^{-1}\times\beta^{-1}\right)  \circ\left(  \alpha
\times\beta\right)  =\operatorname*{id}$. Similarly, $\left(  \alpha
\times\beta\right)  \circ\left(  \alpha^{-1}\times\beta^{-1}\right)
=\operatorname*{id}$. These two equalities show that the maps $\alpha
\times\beta$ and $\alpha^{-1}\times\beta^{-1}$ are mutually inverse. Hence,
the map $\alpha\times\beta:\left[  n\right]  \times\left[  n\right]
\rightarrow\left[  n\right]  \times\left[  n\right]  $ is invertible, and its
inverse is $\alpha^{-1}\times\beta^{-1}$. This proves Lemma
\ref{lem.sol.perm.Inv.sub.bijbij}. \qedhere

\end{vershort}

\begin{verlong}
Now, $\left(  \alpha^{-1}\times\beta^{-1}\right)  \circ\left(  \alpha
\times\beta\right)  =\operatorname*{id}$\ \ \ \ \footnote{\textit{Proof.} Let
$z\in\left[  n\right]  \times\left[  n\right]  $. Thus, $z=\left(  i,j\right)
$ for some $\left(  i,j\right)  \in\left[  n\right]  \times\left[  n\right]
$. Consider this $\left(  i,j\right)  $.
\par
We have%
\begin{align*}
\left(  \left(  \alpha^{-1}\times\beta^{-1}\right)  \circ\left(  \alpha
\times\beta\right)  \right)  \left(  \underbrace{z}_{=\left(  i,j\right)
}\right)   &  =\left(  \left(  \alpha^{-1}\times\beta^{-1}\right)
\circ\left(  \alpha\times\beta\right)  \right)  \left(  i,j\right) \\
&  =\left(  \alpha^{-1}\times\beta^{-1}\right)  \left(  \underbrace{\left(
\alpha\times\beta\right)  \left(  i,j\right)  }_{\substack{=\left(
\alpha\left(  i\right)  ,\beta\left(  j\right)  \right)  \\\text{(by the
definition of }\alpha\times\beta\text{)}}}\right) \\
&  =\left(  \alpha^{-1}\times\beta^{-1}\right)  \left(  \alpha\left(
i\right)  ,\beta\left(  j\right)  \right)  =\left(  \underbrace{\alpha
^{-1}\left(  \alpha\left(  i\right)  \right)  }_{=i},\underbrace{\beta
^{-1}\left(  \beta\left(  j\right)  \right)  }_{=j}\right) \\
&  \ \ \ \ \ \ \ \ \ \ \left(  \text{by the definition of }\alpha^{-1}%
\times\beta^{-1}\right) \\
&  =\left(  i,j\right)  =z=\operatorname*{id}\left(  z\right)  .
\end{align*}
\par
Now, forget that we fixed $z$. We thus have shown that $\left(  \left(
\alpha^{-1}\times\beta^{-1}\right)  \circ\left(  \alpha\times\beta\right)
\right)  \left(  z\right)  =\operatorname*{id}\left(  z\right)  $ for every
$z\in\left[  n\right]  \times\left[  n\right]  $. In other words, $\left(
\alpha^{-1}\times\beta^{-1}\right)  \circ\left(  \alpha\times\beta\right)
=\operatorname*{id}$, qed.} and $\left(  \alpha\times\beta\right)
\circ\left(  \alpha^{-1}\times\beta^{-1}\right)  =\operatorname*{id}%
$\ \ \ \ \footnote{\textit{Proof.} Let $z\in\left[  n\right]  \times\left[
n\right]  $. Thus, $z=\left(  i,j\right)  $ for some $\left(  i,j\right)
\in\left[  n\right]  \times\left[  n\right]  $. Consider this $\left(
i,j\right)  $.
\par
We have%
\begin{align*}
\left(  \left(  \alpha\times\beta\right)  \circ\left(  \alpha^{-1}\times
\beta^{-1}\right)  \right)  \left(  \underbrace{z}_{=\left(  i,j\right)
}\right)   &  =\left(  \left(  \alpha\times\beta\right)  \circ\left(
\alpha^{-1}\times\beta^{-1}\right)  \right)  \left(  i,j\right) \\
&  =\left(  \alpha\times\beta\right)  \left(  \underbrace{\left(  \alpha
^{-1}\times\beta^{-1}\right)  \left(  i,j\right)  }_{\substack{=\left(
\alpha^{-1}\left(  i\right)  ,\beta^{-1}\left(  j\right)  \right)  \\\text{(by
the definition of }\alpha^{-1}\times\beta^{-1}\text{)}}}\right) \\
&  =\left(  \alpha\times\beta\right)  \left(  \alpha^{-1}\left(  i\right)
,\beta^{-1}\left(  j\right)  \right)  =\left(  \underbrace{\alpha\left(
\alpha^{-1}\left(  i\right)  \right)  }_{=i},\underbrace{\beta\left(
\beta^{-1}\left(  j\right)  \right)  }_{=j}\right) \\
&  \ \ \ \ \ \ \ \ \ \ \left(  \text{by the definition of }\alpha\times
\beta\right) \\
&  =\left(  i,j\right)  =z=\operatorname*{id}\left(  z\right)  .
\end{align*}
\par
Now, forget that we fixed $z$. We thus have shown that $\left(  \left(
\alpha\times\beta\right)  \circ\left(  \alpha^{-1}\times\beta^{-1}\right)
\right)  \left(  z\right)  =\operatorname*{id}\left(  z\right)  $ for every
$z\in\left[  n\right]  \times\left[  n\right]  $. In other words, $\left(
\alpha\times\beta\right)  \circ\left(  \alpha^{-1}\times\beta^{-1}\right)
=\operatorname*{id}$, qed.}. Thus, the maps $\alpha\times\beta$ and
$\alpha^{-1}\times\beta^{-1}$ are mutually inverse (since $\left(  \alpha
^{-1}\times\beta^{-1}\right)  \circ\left(  \alpha\times\beta\right)
=\operatorname*{id}$ and $\left(  \alpha\times\beta\right)  \circ\left(
\alpha^{-1}\times\beta^{-1}\right)  =\operatorname*{id}$). Hence, the map
$\alpha\times\beta:\left[  n\right]  \times\left[  n\right]  \rightarrow
\left[  n\right]  \times\left[  n\right]  $ is invertible, and its inverse is
$\alpha^{-1}\times\beta^{-1}$. This proves Lemma
\ref{lem.sol.perm.Inv.sub.bijbij}.
\end{verlong}
\end{proof}

\begin{lemma}
\label{lem.sol.perm.Inv.sub.subset}Let $n\in\mathbb{N}$. Let $\sigma\in S_{n}$
and $\tau\in S_{n}$. Then:

\textbf{(a)} We have $\operatorname*{Inv}\left(  \sigma\circ\tau\right)
\setminus\operatorname*{Inv}\tau\subseteq\left(  \tau\times\tau\right)
^{-1}\left(  \operatorname*{Inv}\sigma\right)  $. (Here, $\tau\times\tau$ is
defined as in Definition \ref{def.sol.perm.Inv.sub.aXb}.)

\textbf{(b)} We have $\left\vert \operatorname*{Inv}\left(  \sigma\circ
\tau\right)  \setminus\operatorname*{Inv}\tau\right\vert \leq\left\vert
\operatorname*{Inv}\sigma\right\vert $.
\end{lemma}

\begin{proof}
[Proof of Lemma \ref{lem.sol.perm.Inv.sub.subset}.]Let $\left[  n\right]  $
denote the set $\left\{  1,2,\ldots,n\right\}  $. Lemma
\ref{lem.sol.perm.Inv.sub.bijbij} (applied to $\alpha=\tau$ and $\beta=\tau$)
yields that the map $\tau\times\tau:\left[  n\right]  \times\left[  n\right]
\rightarrow\left[  n\right]  \times\left[  n\right]  $ is invertible, and its
inverse is $\tau^{-1}\times\tau^{-1}$. Thus, $\left(  \tau\times\tau\right)
^{-1}=\tau^{-1}\times\tau^{-1}$.

\textbf{(a)} Let $c\in\operatorname*{Inv}\left(  \sigma\circ\tau\right)
\setminus\operatorname*{Inv}\tau$. Thus, $c\in\operatorname*{Inv}\left(
\sigma\circ\tau\right)  $ but $c\notin\operatorname*{Inv}\tau$.

We have $c\in\operatorname*{Inv}\left(  \sigma\circ\tau\right)  $. In other
words, $c$ is an inversion of $\sigma\circ\tau$ (since $\operatorname*{Inv}%
\left(  \sigma\circ\tau\right)  $ is the set of all inversions of $\sigma
\circ\tau$). In other words, $c$ is a pair $\left(  i,j\right)  $ of integers
satisfying $1\leq i<j\leq n$ and $\left(  \sigma\circ\tau\right)  \left(
i\right)  >\left(  \sigma\circ\tau\right)  \left(  j\right)  $ (by the
definition of an \textquotedblleft inversion of $\sigma\circ\tau
$\textquotedblright). In other words, there exists a pair $\left(  i,j\right)
$ of integers satisfying $1\leq i<j\leq n$, $\left(  \sigma\circ\tau\right)
\left(  i\right)  >\left(  \sigma\circ\tau\right)  \left(  j\right)  $ and
$c=\left(  i,j\right)  $. Let us denote this pair $\left(  i,j\right)  $ by
$\left(  u,v\right)  $. Thus, $\left(  u,v\right)  $ is a pair of integers
satisfying $1\leq u<v\leq n$, $\left(  \sigma\circ\tau\right)  \left(
u\right)  >\left(  \sigma\circ\tau\right)  \left(  v\right)  $ and $c=\left(
u,v\right)  $.

But $\tau\left(  u\right)  <\tau\left(  v\right)  $%
\ \ \ \ \footnote{\textit{Proof.} Assume the contrary. Thus, $\tau\left(
u\right)  \geq\tau\left(  v\right)  $. But $\tau\in S_{n}$. In other words,
$\tau$ is a permutation of the set $\left\{  1,2,\ldots,n\right\}  $ (since
$S_{n}$ is the set of all permutations of the set $\left\{  1,2,\ldots
,n\right\}  $). Thus, the map $\tau$ is bijective, and therefore also
injective. But $u<v$, so that $u\neq v$ and therefore $\tau\left(  u\right)
\neq\tau\left(  v\right)  $ (since $\tau$ is injective). Combining this with
$\tau\left(  u\right)  \geq\tau\left(  v\right)  $, we obtain $\tau\left(
u\right)  >\tau\left(  v\right)  $.
\par
Now, we know that $\left(  u,v\right)  $ is a pair of integers and satisfies
$1\leq u<v\leq n$ and $\tau\left(  u\right)  >\tau\left(  v\right)  $. In
other words, $\left(  u,v\right)  $ is a pair $\left(  i,j\right)  $ of
integers satisfying $1\leq i<j\leq n$ and $\tau\left(  i\right)  >\tau\left(
j\right)  $. In other words, $\left(  u,v\right)  $ is an inversion of $\tau$
(by the definition of an \textquotedblleft inversion of $\tau$%
\textquotedblright). In other words, $\left(  u,v\right)  \in
\operatorname*{Inv}\tau$ (since $\operatorname*{Inv}\tau$ is the set of all
inversions of $\tau$). But this contradicts $\left(  u,v\right)
=c\notin\operatorname*{Inv}\tau$. This contradiction shows that our assumption
was false, qed.}. Also, $1\leq\tau\left(  u\right)  $ (since $\tau\left(
u\right)  \in\left\{  1,2,\ldots,n\right\}  $) and $\tau\left(  v\right)  \leq
n$ (since $\tau\left(  v\right)  \in\left\{  1,2,\ldots,n\right\}  $).
Finally, $\sigma\left(  \tau\left(  u\right)  \right)  =\left(  \sigma
\circ\tau\right)  \left(  u\right)  >\left(  \sigma\circ\tau\right)  \left(
v\right)  =\sigma\left(  \tau\left(  v\right)  \right)  $.

Thus, $\left(  \tau\left(  u\right)  ,\tau\left(  v\right)  \right)  $ is a
pair of integers satisfying $1\leq\tau\left(  u\right)  <\tau\left(  v\right)
\leq n$ and $\sigma\left(  \tau\left(  u\right)  \right)  >\sigma\left(
\tau\left(  v\right)  \right)  $. In other words, $\left(  \tau\left(
u\right)  ,\tau\left(  v\right)  \right)  $ is a pair $\left(  i,j\right)  $
of integers satisfying $1\leq i<j\leq n$ and $\sigma\left(  i\right)
>\sigma\left(  j\right)  $. In other words, $\left(  \tau\left(  u\right)
,\tau\left(  v\right)  \right)  $ is an inversion of $\sigma$ (by the
definition of an \textquotedblleft inversion of $\sigma$\textquotedblright).
In other words, $\left(  \tau\left(  u\right)  ,\tau\left(  v\right)  \right)
\in\operatorname*{Inv}\sigma$ (since $\operatorname*{Inv}\sigma$ is the set of
all inversions of $\sigma$).

Now, $c=\left(  u,v\right)  \in\left[  n\right]  \times\left[  n\right]  $
(since both $u$ and $v$ belong to $\left[  n\right]  $), and we have
\begin{align*}
\left(  \tau\times\tau\right)  \left(  \underbrace{c}_{=\left(  u,v\right)
}\right)   &  =\left(  \tau\times\tau\right)  \left(  u,v\right)  =\left(
\tau\left(  u\right)  ,\tau\left(  v\right)  \right)
\ \ \ \ \ \ \ \ \ \ \left(  \text{by the definition of }\tau\times\tau\right)
\\
&  \in\operatorname*{Inv}\sigma,
\end{align*}
so that $c\in\left(  \tau\times\tau\right)  ^{-1}\left(  \operatorname*{Inv}%
\sigma\right)  $.

Now, forget that we fixed $c$. We thus have shown that $c\in\left(  \tau
\times\tau\right)  ^{-1}\left(  \operatorname*{Inv}\sigma\right)  $ for every
$c\in\operatorname*{Inv}\left(  \sigma\circ\tau\right)  \setminus
\operatorname*{Inv}\tau$. In other words, $\operatorname*{Inv}\left(
\sigma\circ\tau\right)  \setminus\operatorname*{Inv}\tau\subseteq\left(
\tau\times\tau\right)  ^{-1}\left(  \operatorname*{Inv}\sigma\right)  $. This
proves Lemma \ref{lem.sol.perm.Inv.sub.subset} \textbf{(a)}.

\textbf{(b)} Lemma \ref{lem.sol.perm.Inv.sub.subset} \textbf{(a)} yields
$\operatorname*{Inv}\left(  \sigma\circ\tau\right)  \setminus
\operatorname*{Inv}\tau\subseteq\left(  \tau\times\tau\right)  ^{-1}\left(
\operatorname*{Inv}\sigma\right)  $. Thus,%
\[
\left\vert \operatorname*{Inv}\left(  \sigma\circ\tau\right)  \setminus
\operatorname*{Inv}\tau\right\vert \leq\left\vert \left(  \tau\times
\tau\right)  ^{-1}\left(  \operatorname*{Inv}\sigma\right)  \right\vert
=\left\vert \operatorname*{Inv}\sigma\right\vert
\]
(since $\tau\times\tau$ is a bijection (since the map $\tau\times\tau$ is
invertible)). This proves Lemma \ref{lem.sol.perm.Inv.sub.subset} \textbf{(b)}.
\end{proof}

\begin{lemma}
\label{lem.sol.perm.Inv.sub.l}Let $n\in\mathbb{N}$. Let $\sigma\in S_{n}$.
Then, $\ell\left(  \sigma\right)  =\left\vert \operatorname*{Inv}%
\sigma\right\vert $.
\end{lemma}

\begin{proof}
[Proof of Lemma \ref{lem.sol.perm.Inv.sub.l}.]We have%
\begin{align*}
\ell\left(  \sigma\right)   &  =\left(  \text{the number of inversions of
}\sigma\right)  \ \ \ \ \ \ \ \ \ \ \left(  \text{by the definition of }%
\ell\left(  \sigma\right)  \right) \\
&  =\left\vert \underbrace{\left(  \text{the set of all inversions of }%
\sigma\right)  }_{\substack{=\operatorname*{Inv}\left(  \sigma\right)
\\\text{(since }\operatorname*{Inv}\left(  \sigma\right)  \text{ is the set of
all inversions of }\sigma\text{)}}}\right\vert \\
&  =\left\vert \operatorname*{Inv}\left(  \sigma\right)  \right\vert .
\end{align*}
This proves Lemma \ref{lem.sol.perm.Inv.sub.l}.
\end{proof}

Before we step to the solution to Exercise \ref{exe.perm.Inv.sub}, let us make
a short digression and use Lemma \ref{lem.sol.perm.Inv.sub.subset} to give a
new solution of Exercise \ref{exe.ps2.2.5} \textbf{(c)}:

\begin{proof}
[Third solution to Exercise \ref{exe.ps2.2.5} \textbf{(c)}.]Let $\sigma\in
S_{n}$ and $\tau\in S_{n}$. We have $\left\vert A\setminus B\right\vert
\geq\left\vert A\right\vert -\left\vert B\right\vert $ for any two finite sets
$A$ and $B$\ \ \ \ \footnote{\textit{Proof.} Let $A$ and $B$ be two finite
sets. Then, $A\setminus B=A\setminus\left(  A\cap B\right)  $, so that
\begin{align*}
\left\vert A\setminus B\right\vert  &  =\left\vert A\setminus\left(  A\cap
B\right)  \right\vert =\left\vert A\right\vert -\underbrace{\left\vert A\cap
B\right\vert }_{\substack{\leq\left\vert B\right\vert \\\text{(since }A\cap
B\subseteq B\text{)}}}\ \ \ \ \ \ \ \ \ \ \left(  \text{since }A\cap
B\subseteq A\right) \\
&  \geq\left\vert A\right\vert -\left\vert B\right\vert ,
\end{align*}
qed.}. Applying this to $A=\operatorname*{Inv}\left(  \sigma\circ\tau\right)
$ and $B=\operatorname*{Inv}\tau$, we obtain $\left\vert \operatorname*{Inv}%
\left(  \sigma\circ\tau\right)  \setminus\operatorname*{Inv}\tau\right\vert
\geq\left\vert \operatorname*{Inv}\left(  \sigma\circ\tau\right)  \right\vert
-\left\vert \operatorname*{Inv}\tau\right\vert $. Thus,%
\begin{align*}
\left\vert \operatorname*{Inv}\left(  \sigma\circ\tau\right)  \right\vert
-\left\vert \operatorname*{Inv}\tau\right\vert  &  \leq\left\vert
\operatorname*{Inv}\left(  \sigma\circ\tau\right)  \setminus
\operatorname*{Inv}\tau\right\vert \leq\left\vert \operatorname*{Inv}%
\sigma\right\vert \ \ \ \ \ \ \ \ \ \ \left(  \text{by Lemma
\ref{lem.sol.perm.Inv.sub.subset} \textbf{(b)}}\right) \\
&  =\ell\left(  \sigma\right)  \ \ \ \ \ \ \ \ \ \ \left(  \text{by Lemma
\ref{lem.sol.perm.Inv.sub.l}}\right)  .
\end{align*}
Now,%
\[
\underbrace{\ell\left(  \sigma\circ\tau\right)  }_{\substack{=\left\vert
\operatorname*{Inv}\left(  \sigma\circ\tau\right)  \right\vert \\\text{(by
Lemma \ref{lem.sol.perm.Inv.sub.l}}\\\text{(applied to }\sigma\circ\tau\text{
instead of }\sigma\text{))}}}-\underbrace{\ell\left(  \tau\right)
}_{\substack{=\left\vert \operatorname*{Inv}\tau\right\vert \\\text{(by Lemma
\ref{lem.sol.perm.Inv.sub.l}}\\\text{(applied to }\tau\text{ instead of
}\sigma\text{))}}}=\left\vert \operatorname*{Inv}\left(  \sigma\circ
\tau\right)  \right\vert -\left\vert \operatorname*{Inv}\tau\right\vert
\leq\ell\left(  \sigma\right)  .
\]
In other words, $\ell\left(  \sigma\circ\tau\right)  \leq\ell\left(
\sigma\right)  +\ell\left(  \tau\right)  $. Thus, Exercise \ref{exe.ps2.2.5}
\textbf{(c)} is solved again.
\end{proof}

Now, let us finally solve Exercise \ref{exe.perm.Inv.sub}:

\begin{proof}
[Solution to Exercise \ref{exe.perm.Inv.sub}.]\textbf{(a)} We first observe
that%
\begin{equation}
\underbrace{\ell\left(  \sigma\circ\tau\right)  }_{\substack{=\left\vert
\operatorname*{Inv}\left(  \sigma\circ\tau\right)  \right\vert \\\text{(by
Lemma \ref{lem.sol.perm.Inv.sub.l}}\\\text{(applied to }\sigma\circ
\tau\\\text{instead of }\sigma\text{))}}}-\underbrace{\ell\left(  \tau\right)
}_{\substack{=\left\vert \operatorname*{Inv}\tau\right\vert \\\text{(by Lemma
\ref{lem.sol.perm.Inv.sub.l}}\\\text{(applied to }\tau\\\text{instead of
}\sigma\text{))}}}=\left\vert \operatorname*{Inv}\left(  \sigma\circ
\tau\right)  \right\vert -\left\vert \operatorname*{Inv}\tau\right\vert .
\label{sol.perm.Inv.sub.a.triv}%
\end{equation}

Let us now prove the logical implication%
\begin{equation}
\left(  \ell\left(  \sigma\circ\tau\right)  =\ell\left(  \sigma\right)
+\ell\left(  \tau\right)  \right)  \ \Longrightarrow\ \left(
\operatorname*{Inv}\tau\subseteq\operatorname*{Inv}\left(  \sigma\circ
\tau\right)  \right)  . \label{sol.perm.Inv.sub.a.1}%
\end{equation}

[\textit{Proof of (\ref{sol.perm.Inv.sub.a.1}):} Assume that $\ell\left(
\sigma\circ\tau\right)  =\ell\left(  \sigma\right)  +\ell\left(  \tau\right)
$ holds. We will prove that $\operatorname*{Inv}\tau\subseteq
\operatorname*{Inv}\left(  \sigma\circ\tau\right)  $.

If two finite sets $A$ and $B$ satisfy $\left\vert A\setminus B\right\vert
\leq\left\vert A\right\vert -\left\vert B\right\vert $, then%
\begin{equation}
B\subseteq A \label{sol.perm.Inv.sub.a.1.pf.1}%
\end{equation}
\footnote{\textit{Proof of (\ref{sol.perm.Inv.sub.a.1.pf.1}):} Let $A$ and $B$
be two finite sets satisfying $\left\vert A\setminus B\right\vert
\leq\left\vert A\right\vert -\left\vert B\right\vert $.
\par
We have $A\setminus\left(  A\cap B\right)  =A\setminus B$, so that $\left\vert
A\setminus\left(  A\cap B\right)  \right\vert =\left\vert A\setminus
B\right\vert \leq\left\vert A\right\vert -\left\vert B\right\vert $. Adding
$\left\vert B\right\vert $ to both sides of this inequality, we obtain
$\left\vert A\setminus\left(  A\cap B\right)  \right\vert +\left\vert
B\right\vert \leq\left\vert A\right\vert $. Hence,%
\[
\left\vert A\right\vert \geq\underbrace{\left\vert A\setminus\left(  A\cap
B\right)  \right\vert }_{\substack{=\left\vert A\right\vert -\left\vert A\cap
B\right\vert \\\text{(since }A\cap B\subseteq A\text{)}}}+\left\vert
B\right\vert =\left\vert A\right\vert -\left\vert A\cap B\right\vert
+\left\vert B\right\vert .
\]
Subtracting $\left\vert A\right\vert $ from both sides of this inequality, we
obtain $0\geq-\left\vert A\cap B\right\vert +\left\vert B\right\vert $. In
other words, $\left\vert A\cap B\right\vert \geq\left\vert B\right\vert $.
Also, clearly, $A\cap B$ is a subset of $B$.
\par
But $B$ is a finite set. Hence, the only subset of $B$ having size
$\geq\left\vert B\right\vert $ is $B$ itself. In other words, if $C$ is a
subset of $B$ satisfying $\left\vert C\right\vert \geq\left\vert B\right\vert
$, then $C=B$. Applying this to $C=A\cap B$, we obtain $A\cap B=B$ (since
$A\cap B$ is a subset of $B$ satisfying $\left\vert A\cap B\right\vert
\geq\left\vert B\right\vert $). Hence, $B=A\cap B\subseteq A$. This proves
(\ref{sol.perm.Inv.sub.a.1.pf.1}).}.

Now, Lemma \ref{lem.sol.perm.Inv.sub.subset} \textbf{(b)} yields
\begin{align*}
\left\vert \operatorname*{Inv}\left(  \sigma\circ\tau\right)  \setminus
\operatorname*{Inv}\tau\right\vert  &  \leq\left\vert \operatorname*{Inv}%
\sigma\right\vert =\ell\left(  \sigma\right)  \ \ \ \ \ \ \ \ \ \ \left(
\text{by Lemma \ref{lem.sol.perm.Inv.sub.l}}\right) \\
&  =\ell\left(  \sigma\circ\tau\right)  -\ell\left(  \tau\right)
\ \ \ \ \ \ \ \ \ \ \left(  \text{since }\ell\left(  \sigma\circ\tau\right)
=\ell\left(  \sigma\right)  +\ell\left(  \tau\right)  \right) \\
&  =\left\vert \operatorname*{Inv}\left(  \sigma\circ\tau\right)  \right\vert
-\left\vert \operatorname*{Inv}\tau\right\vert \ \ \ \ \ \ \ \ \ \ \left(
\text{by (\ref{sol.perm.Inv.sub.a.triv})}\right)  .
\end{align*}
Thus, (\ref{sol.perm.Inv.sub.a.1.pf.1}) (applied to $A=\operatorname*{Inv}%
\left(  \sigma\circ\tau\right)  $ and $B=\operatorname*{Inv}\tau$) yields
$\operatorname*{Inv}\tau\subseteq\operatorname*{Inv}\left(  \sigma\circ
\tau\right)  $.

Now, forget our assumption that $\ell\left(  \sigma\circ\tau\right)
=\ell\left(  \sigma\right)  +\ell\left(  \tau\right)  $. We thus have proven
that if $\ell\left(  \sigma\circ\tau\right)  =\ell\left(  \sigma\right)
+\ell\left(  \tau\right)  $, then $\operatorname*{Inv}\tau\subseteq
\operatorname*{Inv}\left(  \sigma\circ\tau\right)  $. In other words, we have
proven the implication (\ref{sol.perm.Inv.sub.a.1}).]

Let us next prove the logical implication%
\begin{equation}
\left(  \operatorname*{Inv}\tau\subseteq\operatorname*{Inv}\left(  \sigma
\circ\tau\right)  \right)  \ \Longrightarrow\ \left(  \ell\left(  \sigma
\circ\tau\right)  =\ell\left(  \sigma\right)  +\ell\left(  \tau\right)
\right)  . \label{sol.perm.Inv.sub.a.2}%
\end{equation}

[\textit{Proof of (\ref{sol.perm.Inv.sub.a.2}):} Assume that
$\operatorname*{Inv}\tau\subseteq\operatorname*{Inv}\left(  \sigma\circ
\tau\right)  $ holds. We will prove that $\ell\left(  \sigma\circ\tau\right)
=\ell\left(  \sigma\right)  +\ell\left(  \tau\right)  $.

Consider the map $\tau\times\tau$ defined as in Definition
\ref{def.sol.perm.Inv.sub.aXb}. Let $\left[  n\right]  $ denote the set
$\left\{  1,2,\ldots,n\right\}  $. Lemma \ref{lem.sol.perm.Inv.sub.bijbij}
(applied to $\alpha=\tau$ and $\beta=\tau$) yields that the map $\tau
\times\tau:\left[  n\right]  \times\left[  n\right]  \rightarrow\left[
n\right]  \times\left[  n\right]  $ is invertible, and its inverse is
$\tau^{-1}\times\tau^{-1}$. Thus, $\left(  \tau\times\tau\right)  ^{-1}%
=\tau^{-1}\times\tau^{-1}$.

Lemma \ref{lem.sol.perm.Inv.sub.subset} \textbf{(a)} shows that%
\begin{equation}
\operatorname*{Inv}\left(  \sigma\circ\tau\right)  \setminus
\operatorname*{Inv}\tau\subseteq\left(  \tau\times\tau\right)  ^{-1}\left(
\operatorname*{Inv}\sigma\right)  . \label{sol.perm.Inv.sub.a.2.pf.1}%
\end{equation}
We shall now prove the reverse inclusion, i.e., we shall prove that $\left(
\tau\times\tau\right)  ^{-1}\left(  \operatorname*{Inv}\sigma\right)
\subseteq\operatorname*{Inv}\left(  \sigma\circ\tau\right)  \setminus
\operatorname*{Inv}\tau$.

Indeed, fix $c\in\left(  \tau\times\tau\right)  ^{-1}\left(
\operatorname*{Inv}\sigma\right)  $. Thus, $c\in\left[  n\right]
\times\left[  n\right]  $ and $\left(  \tau\times\tau\right)  \left(
c\right)  \in\operatorname*{Inv}\sigma$.

We have $\left(  \tau\times\tau\right)  \left(  c\right)  \in
\operatorname*{Inv}\sigma$. In other words, $\left(  \tau\times\tau\right)
\left(  c\right)  $ is an inversion of $\sigma$ (since $\operatorname*{Inv}%
\sigma$ is the set of all inversions of $\sigma$). In other words, $\left(
\tau\times\tau\right)  \left(  c\right)  $ is a pair $\left(  i,j\right)  $ of
integers satisfying $1\leq i<j\leq n$ and $\sigma\left(  i\right)
>\sigma\left(  j\right)  $ (by the definition of an \textquotedblleft
inversion of $\sigma$\textquotedblright). In other words, there exists a pair
$\left(  i,j\right)  $ of integers satisfying $1\leq i<j\leq n$,
$\sigma\left(  i\right)  >\sigma\left(  j\right)  $ and $\left(  \tau
\times\tau\right)  \left(  c\right)  =\left(  i,j\right)  $. Let us denote
this pair $\left(  i,j\right)  $ by $\left(  u,v\right)  $. Thus, $\left(
u,v\right)  $ is a pair of integers satisfying $1\leq u<v\leq n$,
$\sigma\left(  u\right)  >\sigma\left(  v\right)  $ and $\left(  \tau
\times\tau\right)  \left(  c\right)  =\left(  u,v\right)  $.

From $\left(  \tau\times\tau\right)  \left(  c\right)  =\left(  u,v\right)  $,
we obtain%
\begin{align*}
c  &  =\underbrace{\left(  \tau\times\tau\right)  ^{-1}}_{=\tau^{-1}\times
\tau^{-1}}\left(  u,v\right)  \ \ \ \ \ \ \ \ \ \ \left(  \text{since the map
}\tau\times\tau\text{ is invertible}\right) \\
&  =\left(  \tau^{-1}\times\tau^{-1}\right)  \left(  u,v\right)  =\left(
\tau^{-1}\left(  u\right)  ,\tau^{-1}\left(  v\right)  \right)
\end{align*}
(by the definition of $\tau^{-1}\times\tau^{-1}$).

Notice that%
\begin{equation}
\left(  \sigma\circ\tau\right)  \left(  \tau^{-1}\left(  v\right)  \right)
=\sigma\left(  \underbrace{\tau\left(  \tau^{-1}\left(  v\right)  \right)
}_{=v}\right)  =\sigma\left(  v\right)  \label{sol.perm.Inv.sub.a.2.pf.4}%
\end{equation}
and%
\begin{align}
\left(  \sigma\circ\tau\right)  \left(  \tau^{-1}\left(  u\right)  \right)
&  =\sigma\left(  \underbrace{\tau\left(  \tau^{-1}\left(  u\right)  \right)
}_{=u}\right)  =\sigma\left(  u\right) \nonumber\\
&  >\sigma\left(  v\right)  =\left(  \sigma\circ\tau\right)  \left(  \tau
^{-1}\left(  v\right)  \right)  \label{sol.perm.Inv.sub.a.2.pf.5}%
\end{align}
(by (\ref{sol.perm.Inv.sub.a.2.pf.4})).

Let us now prove that $\tau^{-1}\left(  u\right)  <\tau^{-1}\left(  v\right)
$. Indeed, let us assume the contrary. Thus, $\tau^{-1}\left(  u\right)
\geq\tau^{-1}\left(  v\right)  $. But $u\neq v$ (since $u<v$), so that
$\tau^{-1}\left(  u\right)  \neq\tau^{-1}\left(  v\right)  $. Combined with
$\tau^{-1}\left(  u\right)  \geq\tau^{-1}\left(  v\right)  $, this yields
$\tau^{-1}\left(  u\right)  >\tau^{-1}\left(  v\right)  $. In other words,
$\tau^{-1}\left(  v\right)  <\tau^{-1}\left(  u\right)  $. Also, $1\leq
\tau^{-1}\left(  v\right)  $ (since $\tau^{-1}\left(  v\right)  \in\left\{
1,2,\ldots,n\right\}  $) and $\tau^{-1}\left(  u\right)  \leq n$ (since
$\tau^{-1}\left(  u\right)  \in\left\{  1,2,\ldots,n\right\}  $). Finally,
$\tau\left(  \tau^{-1}\left(  v\right)  \right)  =v>u=\tau\left(  \tau
^{-1}\left(  u\right)  \right)  $.

Thus, $\left(  \tau^{-1}\left(  v\right)  ,\tau^{-1}\left(  u\right)  \right)
$ is a pair of integers satisfying $1\leq\tau^{-1}\left(  v\right)  <\tau
^{-1}\left(  u\right)  \leq n$ and $\tau\left(  \tau^{-1}\left(  v\right)
\right)  >\tau\left(  \tau^{-1}\left(  u\right)  \right)  $. In other words,
$\left(  \tau^{-1}\left(  v\right)  ,\tau^{-1}\left(  u\right)  \right)  $ is
a pair $\left(  i,j\right)  $ of integers satisfying $1\leq i<j\leq n$ and
$\tau\left(  i\right)  >\tau\left(  j\right)  $. In other words, $\left(
\tau^{-1}\left(  v\right)  ,\tau^{-1}\left(  u\right)  \right)  $ is an
inversion of $\tau$ (by the definition of an \textquotedblleft inversion of
$\tau$\textquotedblright). In other words, $\left(  \tau^{-1}\left(  v\right)
,\tau^{-1}\left(  u\right)  \right)  \in\operatorname*{Inv}\tau$ (since
$\operatorname*{Inv}\tau$ is the set of all inversions of $\tau$). Hence,
$\left(  \tau^{-1}\left(  v\right)  ,\tau^{-1}\left(  u\right)  \right)
\in\operatorname*{Inv}\tau\subseteq\operatorname*{Inv}\left(  \sigma\circ
\tau\right)  $. In other words, $\left(  \tau^{-1}\left(  v\right)  ,\tau
^{-1}\left(  u\right)  \right)  $ is an inversion of $\sigma\circ\tau$ (since
$\operatorname*{Inv}\left(  \sigma\circ\tau\right)  $ is the set of all
inversions of $\sigma\circ\tau$). In other words, $\left(  \tau^{-1}\left(
v\right)  ,\tau^{-1}\left(  u\right)  \right)  $ is a pair $\left(
i,j\right)  $ of integers satisfying $1\leq i<j\leq n$ and $\left(
\sigma\circ\tau\right)  \left(  i\right)  >\left(  \sigma\circ\tau\right)
\left(  j\right)  $ (by the definition of an \textquotedblleft inversion of
$\sigma\circ\tau$\textquotedblright). In other words, $\left(  \tau
^{-1}\left(  v\right)  ,\tau^{-1}\left(  u\right)  \right)  $ is a pair of
integers satisfying $1\leq\tau^{-1}\left(  v\right)  <\tau^{-1}\left(
u\right)  \leq n$ and $\left(  \sigma\circ\tau\right)  \left(  \tau
^{-1}\left(  v\right)  \right)  >\left(  \sigma\circ\tau\right)  \left(
\tau^{-1}\left(  u\right)  \right)  $. But $\left(  \sigma\circ\tau\right)
\left(  \tau^{-1}\left(  v\right)  \right)  >\left(  \sigma\circ\tau\right)
\left(  \tau^{-1}\left(  u\right)  \right)  >\left(  \sigma\circ\tau\right)
\left(  \tau^{-1}\left(  v\right)  \right)  $ (by
(\ref{sol.perm.Inv.sub.a.2.pf.5})), which is absurd. This contradiction proves
that our assumption was wrong. Hence, $\tau^{-1}\left(  u\right)  <\tau
^{-1}\left(  v\right)  $ is proven.

Now, $1\leq\tau^{-1}\left(  u\right)  $ (since $\tau^{-1}\left(  u\right)
\in\left\{  1,2,\ldots,n\right\}  $) and $\tau^{-1}\left(  v\right)  \leq n$
(since $\tau^{-1}\left(  v\right)  \in\left\{  1,2,\ldots,n\right\}  $).
Finally, recall that (\ref{sol.perm.Inv.sub.a.2.pf.5}) holds. Thus, $\left(
\tau^{-1}\left(  u\right)  ,\tau^{-1}\left(  v\right)  \right)  $ is a pair of
integers satisfying $1\leq\tau^{-1}\left(  u\right)  <\tau^{-1}\left(
v\right)  \leq n$ and $\left(  \sigma\circ\tau\right)  \left(  \tau
^{-1}\left(  u\right)  \right)  >\left(  \sigma\circ\tau\right)  \left(
\tau^{-1}\left(  v\right)  \right)  $. In other words, $\left(  \tau
^{-1}\left(  u\right)  ,\tau^{-1}\left(  v\right)  \right)  $ is a pair
$\left(  i,j\right)  $ of integers satisfying $1\leq i<j\leq n$ and $\left(
\sigma\circ\tau\right)  \left(  i\right)  >\left(  \sigma\circ\tau\right)
\left(  j\right)  $. In other words, $\left(  \tau^{-1}\left(  u\right)
,\tau^{-1}\left(  v\right)  \right)  $ is an inversion of $\sigma\circ\tau$
(by the definition of an \textquotedblleft inversion of $\sigma\circ\tau
$\textquotedblright). In other words, $\left(  \tau^{-1}\left(  u\right)
,\tau^{-1}\left(  v\right)  \right)  \in\operatorname*{Inv}\left(  \sigma
\circ\tau\right)  $ (since $\operatorname*{Inv}\left(  \sigma\circ\tau\right)
$ is the set of all inversions of $\sigma\circ\tau$).

Hence, $c=\left(  \tau^{-1}\left(  u\right)  ,\tau^{-1}\left(  v\right)
\right)  \in\operatorname*{Inv}\left(  \sigma\circ\tau\right)  $.

Let us now prove that $c\notin\operatorname*{Inv}\tau$. Indeed, assume the
contrary. Thus, $c\in\operatorname*{Inv}\tau$. Thus, $\left(  \tau^{-1}\left(
u\right)  ,\tau^{-1}\left(  v\right)  \right)  =c\in\operatorname*{Inv}\tau$.
In other words, $\left(  \tau^{-1}\left(  u\right)  ,\tau^{-1}\left(
v\right)  \right)  $ is an inversion of $\tau$ (since $\operatorname*{Inv}%
\tau$ is the set of all inversions of $\tau$). In other words, $\left(
\tau^{-1}\left(  u\right)  ,\tau^{-1}\left(  v\right)  \right)  $ is a pair
$\left(  i,j\right)  $ of integers satisfying $1\leq i<j\leq n$ and
$\tau\left(  i\right)  >\tau\left(  j\right)  $ (by the definition of an
\textquotedblleft inversion of $\tau$\textquotedblright). In other words,
$\left(  \tau^{-1}\left(  u\right)  ,\tau^{-1}\left(  v\right)  \right)  $ is
a pair of integers satisfying $1\leq\tau^{-1}\left(  u\right)  <\tau
^{-1}\left(  v\right)  \leq n$ and $\tau\left(  \tau^{-1}\left(  u\right)
\right)  >\tau\left(  \tau^{-1}\left(  v\right)  \right)  $. But $\tau\left(
\tau^{-1}\left(  u\right)  \right)  >\tau\left(  \tau^{-1}\left(  v\right)
\right)  =v$ contradicts $\tau\left(  \tau^{-1}\left(  u\right)  \right)
=u<v$. This contradiction proves that our assumption was wrong. Hence,
$c\notin\operatorname*{Inv}\tau$ is proven.

Combining $c\in\operatorname*{Inv}\left(  \sigma\circ\tau\right)  $ with
$c\notin\operatorname*{Inv}\tau$, we obtain $c\in\operatorname*{Inv}\left(
\sigma\circ\tau\right)  \setminus\operatorname*{Inv}\tau$.

Now, forget that we fixed $c$. We thus have shown that $c\in
\operatorname*{Inv}\left(  \sigma\circ\tau\right)  \setminus
\operatorname*{Inv}\tau$ for each $c\in\left(  \tau\times\tau\right)
^{-1}\left(  \operatorname*{Inv}\sigma\right)  $. In other words, we have
proven the inclusion%
\[
\left(  \tau\times\tau\right)  ^{-1}\left(  \operatorname*{Inv}\sigma\right)
\subseteq\operatorname*{Inv}\left(  \sigma\circ\tau\right)  \setminus
\operatorname*{Inv}\tau.
\]
Combining this with the inclusion (\ref{sol.perm.Inv.sub.a.2.pf.1}), we obtain%
\[
\operatorname*{Inv}\left(  \sigma\circ\tau\right)  \setminus
\operatorname*{Inv}\tau=\left(  \tau\times\tau\right)  ^{-1}\left(
\operatorname*{Inv}\sigma\right)  .
\]
Hence,%
\[
\left\vert \operatorname*{Inv}\left(  \sigma\circ\tau\right)  \setminus
\operatorname*{Inv}\tau\right\vert =\left\vert \left(  \tau\times\tau\right)
^{-1}\left(  \operatorname*{Inv}\sigma\right)  \right\vert =\left\vert
\operatorname*{Inv}\sigma\right\vert
\]
(since the map $\tau\times\tau$ is a bijection (since $\tau\times\tau$ is
invertible)). Thus,
\[
\left\vert \operatorname*{Inv}\left(  \sigma\circ\tau\right)  \setminus
\operatorname*{Inv}\tau\right\vert =\left\vert \operatorname*{Inv}%
\sigma\right\vert =\ell\left(  \sigma\right)
\]
(by Lemma \ref{lem.sol.perm.Inv.sub.l}). Hence,%
\begin{align*}
\ell\left(  \sigma\right)   &  =\left\vert \operatorname*{Inv}\left(
\sigma\circ\tau\right)  \setminus\operatorname*{Inv}\tau\right\vert
=\left\vert \operatorname*{Inv}\left(  \sigma\circ\tau\right)  \right\vert
-\left\vert \operatorname*{Inv}\tau\right\vert \ \ \ \ \ \ \ \ \ \ \left(
\text{since }\operatorname*{Inv}\tau\subseteq\operatorname*{Inv}\left(
\sigma\circ\tau\right)  \right) \\
&  =\ell\left(  \sigma\circ\tau\right)  -\ell\left(  \tau\right)
\ \ \ \ \ \ \ \ \ \ \left(  \text{by (\ref{sol.perm.Inv.sub.a.triv})}\right)
.
\end{align*}
In other words, $\ell\left(  \sigma\circ\tau\right)  =\ell\left(
\sigma\right)  +\ell\left(  \tau\right)  $.

Now, forget our assumption that $\operatorname*{Inv}\tau\subseteq
\operatorname*{Inv}\left(  \sigma\circ\tau\right)  $. We thus have proven that
if $\operatorname*{Inv}\tau\subseteq\operatorname*{Inv}\left(  \sigma\circ
\tau\right)  $, then $\ell\left(  \sigma\circ\tau\right)  =\ell\left(
\sigma\right)  +\ell\left(  \tau\right)  $. In other words, we have proven the
implication (\ref{sol.perm.Inv.sub.a.2}).]

Combining the two implications (\ref{sol.perm.Inv.sub.a.1}) and
(\ref{sol.perm.Inv.sub.a.2}), we obtain the logical equivalence%
\[
\left(  \ell\left(  \sigma\circ\tau\right)  =\ell\left(  \sigma\right)
+\ell\left(  \tau\right)  \right)  \ \Longleftrightarrow\ \left(
\operatorname*{Inv}\tau\subseteq\operatorname*{Inv}\left(  \sigma\circ
\tau\right)  \right)  .
\]
In other words, $\ell\left(  \sigma\circ\tau\right)  =\ell\left(
\sigma\right)  +\ell\left(  \tau\right)  $ holds if and only if
$\operatorname*{Inv}\tau\subseteq\operatorname*{Inv}\left(  \sigma\circ
\tau\right)  $. This solves Exercise \ref{exe.perm.Inv.sub} \textbf{(a)}.

\textbf{(b)} Exercise \ref{exe.ps2.2.5} \textbf{(f)} (applied to $\sigma
\circ\tau$ instead of $\sigma$) yields $\ell\left(  \sigma\circ\tau\right)
=\ell\left(  \underbrace{\left(  \sigma\circ\tau\right)  ^{-1}}_{=\tau
^{-1}\circ\sigma^{-1}}\right)  =\ell\left(  \tau^{-1}\circ\sigma^{-1}\right)
$. Exercise \ref{exe.ps2.2.5} \textbf{(f)} (applied to $\tau$ instead of
$\sigma$) yields $\ell\left(  \tau\right)  =\ell\left(  \tau^{-1}\right)  $.
Exercise \ref{exe.ps2.2.5} \textbf{(f)} yields $\ell\left(  \sigma\right)
=\ell\left(  \sigma^{-1}\right)  $. Hence, $\underbrace{\ell\left(
\sigma\right)  }_{=\ell\left(  \sigma^{-1}\right)  }+\underbrace{\ell\left(
\tau\right)  }_{=\ell\left(  \tau^{-1}\right)  }=\ell\left(  \sigma
^{-1}\right)  +\ell\left(  \tau^{-1}\right)  =\ell\left(  \tau^{-1}\right)
+\ell\left(  \sigma^{-1}\right)  $.

Exercise \ref{exe.perm.Inv.sub} \textbf{(a)} (applied to $\tau^{-1}$ and
$\sigma^{-1}$ instead of $\sigma$ and $\tau$) yields that $\ell\left(
\tau^{-1}\circ\sigma^{-1}\right)  =\ell\left(  \tau^{-1}\right)  +\ell\left(
\sigma^{-1}\right)  $ holds if and only if $\operatorname*{Inv}\left(
\sigma^{-1}\right)  \subseteq\operatorname*{Inv}\left(  \tau^{-1}\circ
\sigma^{-1}\right)  $. In light of the equalities $\ell\left(  \tau^{-1}%
\circ\sigma^{-1}\right)  =\ell\left(  \sigma\circ\tau\right)  $ and
$\ell\left(  \tau^{-1}\right)  +\ell\left(  \sigma^{-1}\right)  =\ell\left(
\sigma\right)  +\ell\left(  \tau\right)  $, we can rewrite this as follows:%
\[
\ell\left(  \sigma\circ\tau\right)  =\ell\left(  \sigma\right)  +\ell\left(
\tau\right)  \text{ holds if and only if }\operatorname*{Inv}\left(
\sigma^{-1}\right)  \subseteq\operatorname*{Inv}\left(  \tau^{-1}\circ
\sigma^{-1}\right)  .
\]
This solves Exercise \ref{exe.perm.Inv.sub} \textbf{(b)}.

\textbf{(c)} Exercise \ref{exe.perm.Inv.sub} \textbf{(a)} (applied to
$\tau\circ\sigma^{-1}$ and $\sigma$ instead of $\sigma$ and $\tau$) yields
that $\ell\left(  \tau\circ\sigma^{-1}\circ\sigma\right)  =\ell\left(
\tau\circ\sigma^{-1}\right)  +\ell\left(  \sigma\right)  $ holds if and only
if $\operatorname*{Inv}\sigma\subseteq\operatorname*{Inv}\left(  \tau
\circ\sigma^{-1}\circ\sigma\right)  $. In light of $\tau\circ
\underbrace{\sigma^{-1}\circ\sigma}_{=\operatorname*{id}}=\tau\circ
\operatorname*{id}=\tau$, this rewrites as follows:%
\[
\ell\left(  \tau\right)  =\ell\left(  \tau\circ\sigma^{-1}\right)
+\ell\left(  \sigma\right)  \text{ holds if and only if }\operatorname*{Inv}%
\sigma\subseteq\operatorname*{Inv}\tau.
\]
In other words, $\operatorname*{Inv}\sigma\subseteq\operatorname*{Inv}\tau$
holds if and only if $\ell\left(  \tau\right)  =\ell\left(  \tau\circ
\sigma^{-1}\right)  +\ell\left(  \sigma\right)  $. This solves Exercise
\ref{exe.perm.Inv.sub} \textbf{(c)}.

\textbf{(d)} Assume that $\operatorname*{Inv}\sigma=\operatorname*{Inv}\tau$.
We WLOG assume that $\ell\left(  \sigma\right)  \geq\ell\left(  \tau\right)  $
(since otherwise, we can simply swap $\sigma$ with $\tau$). We have
$\operatorname*{Inv}\sigma\subseteq\operatorname*{Inv}\tau$ (since
$\operatorname*{Inv}\sigma=\operatorname*{Inv}\tau$). But Exercise
\ref{exe.perm.Inv.sub} \textbf{(c)} shows that $\operatorname*{Inv}%
\sigma\subseteq\operatorname*{Inv}\tau$ holds if and only if $\ell\left(
\tau\right)  =\ell\left(  \tau\circ\sigma^{-1}\right)  +\ell\left(
\sigma\right)  $. Hence, we have $\ell\left(  \tau\right)  =\ell\left(
\tau\circ\sigma^{-1}\right)  +\ell\left(  \sigma\right)  $ (since
$\operatorname*{Inv}\sigma\subseteq\operatorname*{Inv}\tau$ holds). Thus,
$\ell\left(  \tau\right)  =\ell\left(  \tau\circ\sigma^{-1}\right)
+\underbrace{\ell\left(  \sigma\right)  }_{\geq\ell\left(  \tau\right)  }%
\geq\ell\left(  \tau\circ\sigma^{-1}\right)  +\ell\left(  \tau\right)  $.
Subtracting $\ell\left(  \tau\right)  $ from both sides of this inequality, we
obtain $0\geq\ell\left(  \tau\circ\sigma^{-1}\right)  $. In other words,
$\ell\left(  \tau\circ\sigma^{-1}\right)  \leq0$.

But $\ell\left(  \tau\circ\sigma^{-1}\right)  $ is the number of inversions of
$\tau\circ\sigma^{-1}$ (by the definition of $\ell\left(  \tau\circ\sigma
^{-1}\right)  $), and thus is a nonnegative integer. Hence, $\ell\left(
\tau\circ\sigma^{-1}\right)  \geq0$. Combining this with $\ell\left(
\tau\circ\sigma^{-1}\right)  \leq0$, we obtain $\ell\left(  \tau\circ
\sigma^{-1}\right)  =0$. Thus, Corollary \ref{cor.sol.ps2.2.5.d2} (applied to
$\tau\circ\sigma^{-1}$ instead of $\sigma$) yields $\tau\circ\sigma
^{-1}=\operatorname*{id}$. Hence, $\underbrace{\tau\circ\sigma^{-1}%
}_{=\operatorname*{id}}\circ\sigma=\operatorname*{id}\circ\sigma=\sigma$, so
that $\sigma=\tau\circ\underbrace{\sigma^{-1}\circ\sigma}_{=\operatorname*{id}%
}=\tau\circ\operatorname*{id}=\tau$. This solves Exercise
\ref{exe.perm.Inv.sub} \textbf{(d)}.

\textbf{(e)} Let us first prove the logical implication%
\begin{equation}
\left(  \operatorname*{Inv}\tau\subseteq\operatorname*{Inv}\left(  \sigma
\circ\tau\right)  \right)  \ \Longrightarrow\ \left(  \left(
\operatorname*{Inv}\sigma\right)  \cap\left(  \operatorname*{Inv}\left(
\tau^{-1}\right)  \right)  =\varnothing\right)  . \label{sol.perm.Inv.sub.e.1}%
\end{equation}

[\textit{Proof of (\ref{sol.perm.Inv.sub.e.1}):} Assume that
$\operatorname*{Inv}\tau\subseteq\operatorname*{Inv}\left(  \sigma\circ
\tau\right)  $ holds. We will prove that $\left(  \operatorname*{Inv}%
\sigma\right)  \cap\left(  \operatorname*{Inv}\left(  \tau^{-1}\right)
\right)  =\varnothing$.

Indeed, fix $c\in\left(  \operatorname*{Inv}\sigma\right)  \cap\left(
\operatorname*{Inv}\left(  \tau^{-1}\right)  \right)  $. Thus, $c\in\left(
\operatorname*{Inv}\sigma\right)  \cap\left(  \operatorname*{Inv}\left(
\tau^{-1}\right)  \right)  \subseteq\operatorname*{Inv}\sigma$ and
$c\in\left(  \operatorname*{Inv}\sigma\right)  \cap\left(  \operatorname*{Inv}%
\left(  \tau^{-1}\right)  \right)  \subseteq\operatorname*{Inv}\left(
\tau^{-1}\right)  $.

We have $c\in\operatorname*{Inv}\sigma$. In other words, $c$ is an inversion
of $\sigma$ (since $\operatorname*{Inv}\sigma$ is the set of all inversions of
$\sigma$). In other words, $c$ is a pair $\left(  i,j\right)  $ of integers
satisfying $1\leq i<j\leq n$ and $\sigma\left(  i\right)  >\sigma\left(
j\right)  $ (by the definition of an \textquotedblleft inversion of $\sigma
$\textquotedblright). In other words, there exists a pair $\left(  i,j\right)
$ of integers satisfying $1\leq i<j\leq n$, $\sigma\left(  i\right)
>\sigma\left(  j\right)  $ and $c=\left(  i,j\right)  $. Let us denote this
pair $\left(  i,j\right)  $ by $\left(  u,v\right)  $. Thus, $\left(
u,v\right)  $ is a pair of integers satisfying $1\leq u<v\leq n$,
$\sigma\left(  u\right)  >\sigma\left(  v\right)  $ and $c=\left(  u,v\right)
$. We have $v>u$ (since $u<v$). Thus, $\tau\left(  \tau^{-1}\left(  v\right)
\right)  =v>u=\tau\left(  \tau^{-1}\left(  u\right)  \right)  $.

We have $\left(  u,v\right)  =c\in\operatorname*{Inv}\left(  \tau^{-1}\right)
$. In other words, $\left(  u,v\right)  $ is an inversion of $\tau^{-1}$
(since $\operatorname*{Inv}\left(  \tau^{-1}\right)  $ is the set of all
inversions of $\tau^{-1}$). In other words, $\left(  u,v\right)  $ is a pair
$\left(  i,j\right)  $ of integers satisfying $1\leq i<j\leq n$ and $\tau
^{-1}\left(  i\right)  >\tau^{-1}\left(  j\right)  $ (by the definition of an
\textquotedblleft inversion of $\tau^{-1}$\textquotedblright). In other words,
$\left(  u,v\right)  $ is a pair of integers satisfying $1\leq u<v\leq n$ and
$\tau^{-1}\left(  u\right)  >\tau^{-1}\left(  v\right)  $. From $\tau
^{-1}\left(  u\right)  >\tau^{-1}\left(  v\right)  $, we obtain $\tau
^{-1}\left(  v\right)  <\tau^{-1}\left(  u\right)  $.

We have $\tau^{-1}\left(  u\right)  \in\left\{  1,2,\ldots,n\right\}  $ (since
$\tau\in S_{n}$), so that $\tau^{-1}\left(  u\right)  \leq n$. Also,
$\tau^{-1}\left(  v\right)  \in\left\{  1,2,\ldots,n\right\}  $ (since
$\tau\in S_{n}$) and thus $1\leq\tau^{-1}\left(  v\right)  $. Altogether, we
thus know that $\left(  \tau^{-1}\left(  v\right)  ,\tau^{-1}\left(  u\right)
\right)  $ is a pair of integers satisfying $1\leq\tau^{-1}\left(  v\right)
<\tau^{-1}\left(  u\right)  \leq n$ and $\tau\left(  \tau^{-1}\left(
v\right)  \right)  >\tau\left(  \tau^{-1}\left(  u\right)  \right)  $. In
other words, $\left(  \tau^{-1}\left(  v\right)  ,\tau^{-1}\left(  u\right)
\right)  $ is a pair $\left(  i,j\right)  $ of integers satisfying $1\leq
i<j\leq n$ and $\tau\left(  i\right)  >\tau\left(  j\right)  $. In other
words, $\left(  \tau^{-1}\left(  v\right)  ,\tau^{-1}\left(  u\right)
\right)  $ is an inversion of $\tau$ (by the definition of an
\textquotedblleft inversion of $\tau$\textquotedblright). In other words,
$\left(  \tau^{-1}\left(  v\right)  ,\tau^{-1}\left(  u\right)  \right)
\in\operatorname*{Inv}\tau$ (since $\operatorname*{Inv}\tau$ is the set of all
inversions of $\tau$).

Hence, $\left(  \tau^{-1}\left(  v\right)  ,\tau^{-1}\left(  u\right)
\right)  \in\operatorname*{Inv}\tau\subseteq\operatorname*{Inv}\left(
\sigma\circ\tau\right)  $. In other words, $\left(  \tau^{-1}\left(  v\right)
,\tau^{-1}\left(  u\right)  \right)  $ is an inversion of $\sigma\circ\tau$
(since $\operatorname*{Inv}\left(  \sigma\circ\tau\right)  $ is the set of all
inversions of $\sigma\circ\tau$). In other words, $\left(  \tau^{-1}\left(
v\right)  ,\tau^{-1}\left(  u\right)  \right)  $ is a pair $\left(
i,j\right)  $ of integers satisfying $1\leq i<j\leq n$ and $\left(
\sigma\circ\tau\right)  \left(  i\right)  >\left(  \sigma\circ\tau\right)
\left(  j\right)  $ (by the definition of an \textquotedblleft inversion of
$\sigma\circ\tau$\textquotedblright). In other words, $\left(  \tau
^{-1}\left(  v\right)  ,\tau^{-1}\left(  u\right)  \right)  $ is a pair of
integers satisfying $1\leq\tau^{-1}\left(  v\right)  <\tau^{-1}\left(
u\right)  \leq n$ and $\left(  \sigma\circ\tau\right)  \left(  \tau
^{-1}\left(  v\right)  \right)  >\left(  \sigma\circ\tau\right)  \left(
\tau^{-1}\left(  u\right)  \right)  $.

In particular, we have $\left(  \sigma\circ\tau\right)  \left(  \tau
^{-1}\left(  v\right)  \right)  >\left(  \sigma\circ\tau\right)  \left(
\tau^{-1}\left(  u\right)  \right)  $. This rewrites as $\sigma\left(
v\right)  >\sigma\left(  u\right)  $ (since $\left(  \sigma\circ\tau\right)
\left(  \tau^{-1}\left(  v\right)  \right)  =\left(  \sigma\circ
\underbrace{\tau\circ\tau^{-1}}_{=\operatorname*{id}}\right)  \left(
v\right)  =\sigma\left(  v\right)  $ and $\left(  \sigma\circ\tau\right)
\left(  \tau^{-1}\left(  u\right)  \right)  =\left(  \sigma\circ
\underbrace{\tau\circ\tau^{-1}}_{=\operatorname*{id}}\right)  \left(
u\right)  =\sigma\left(  u\right)  $). But this contradicts $\sigma\left(
u\right)  >\sigma\left(  v\right)  $.

Now, forget that we fixed $c$. We thus have found a contradiction for each
$c\in\left(  \operatorname*{Inv}\sigma\right)  \cap\left(  \operatorname*{Inv}%
\left(  \tau^{-1}\right)  \right)  $. Hence, there exists no $c\in\left(
\operatorname*{Inv}\sigma\right)  \cap\left(  \operatorname*{Inv}\left(
\tau^{-1}\right)  \right)  $. In other words, $\left(  \operatorname*{Inv}%
\sigma\right)  \cap\left(  \operatorname*{Inv}\left(  \tau^{-1}\right)
\right)  $ is the empty set. In other words, $\left(  \operatorname*{Inv}%
\sigma\right)  \cap\left(  \operatorname*{Inv}\left(  \tau^{-1}\right)
\right)  =\varnothing$.

Now, forget our assumption that $\operatorname*{Inv}\tau\subseteq
\operatorname*{Inv}\left(  \sigma\circ\tau\right)  $. We thus have proven that
if $\operatorname*{Inv}\tau\subseteq\operatorname*{Inv}\left(  \sigma\circ
\tau\right)  $, then $\left(  \operatorname*{Inv}\sigma\right)  \cap\left(
\operatorname*{Inv}\left(  \tau^{-1}\right)  \right)  =\varnothing$. In other
words, we have proven the implication (\ref{sol.perm.Inv.sub.e.1}).]

Next, we shall prove the logical implication%
\begin{equation}
\left(  \left(  \operatorname*{Inv}\sigma\right)  \cap\left(
\operatorname*{Inv}\left(  \tau^{-1}\right)  \right)  =\varnothing\right)
\ \Longrightarrow\ \left(  \operatorname*{Inv}\tau\subseteq\operatorname*{Inv}%
\left(  \sigma\circ\tau\right)  \right)  . \label{sol.perm.Inv.sub.e.2}%
\end{equation}

[\textit{Proof of (\ref{sol.perm.Inv.sub.e.2}):} Assume that $\left(
\operatorname*{Inv}\sigma\right)  \cap\left(  \operatorname*{Inv}\left(
\tau^{-1}\right)  \right)  =\varnothing$ holds. We will prove that
$\operatorname*{Inv}\tau\subseteq\operatorname*{Inv}\left(  \sigma\circ
\tau\right)  $.

Indeed, let $c\in\left(  \operatorname*{Inv}\tau\right)  \setminus\left(
\operatorname*{Inv}\left(  \sigma\circ\tau\right)  \right)  $. We shall prove
a contradiction.

We have $c\in\left(  \operatorname*{Inv}\tau\right)  \setminus\left(
\operatorname*{Inv}\left(  \sigma\circ\tau\right)  \right)  \subseteq
\operatorname*{Inv}\tau$. In other words, $c$ is an inversion of $\tau$ (since
$\operatorname*{Inv}\tau$ is the set of all inversions of $\tau$). In other
words, $c$ is a pair $\left(  i,j\right)  $ of integers satisfying $1\leq
i<j\leq n$ and $\tau\left(  i\right)  >\tau\left(  j\right)  $ (by the
definition of an \textquotedblleft inversion of $\tau$\textquotedblright). In
other words, there exists a pair $\left(  i,j\right)  $ of integers satisfying
$1\leq i<j\leq n$, $\tau\left(  i\right)  >\tau\left(  j\right)  $ and
$c=\left(  i,j\right)  $. Let us denote this pair $\left(  i,j\right)  $ by
$\left(  u,v\right)  $. Thus, $\left(  u,v\right)  $ is a pair of integers
satisfying $1\leq u<v\leq n$, $\tau\left(  u\right)  >\tau\left(  v\right)  $
and $c=\left(  u,v\right)  $. We have $v>u$ (since $u<v$) and $\tau\left(
v\right)  <\tau\left(  u\right)  $ (since $\tau\left(  u\right)  >\tau\left(
v\right)  $).

From $\tau\in S_{n}$, we obtain $\tau\left(  v\right)  \in\left\{
1,2,\ldots,n\right\}  $ and thus $1\leq\tau\left(  v\right)  $. From $\tau\in
S_{n}$, we obtain $\tau\left(  u\right)  \in\left\{  1,2,\ldots,n\right\}  $
and thus $\tau\left(  u\right)  \leq n$. Also, $\tau^{-1}\left(  \tau\left(
v\right)  \right)  =v>u=\tau^{-1}\left(  \tau\left(  u\right)  \right)  $.
Altogether, we thus know that $\left(  \tau\left(  v\right)  ,\tau\left(
u\right)  \right)  $ is a pair of integers satisfying $1\leq\tau\left(
v\right)  <\tau\left(  u\right)  \leq n$ and $\tau^{-1}\left(  \tau\left(
v\right)  \right)  >\tau^{-1}\left(  \tau\left(  u\right)  \right)  $. In
other words, $\left(  \tau\left(  v\right)  ,\tau\left(  u\right)  \right)  $
is a pair $\left(  i,j\right)  $ of integers satisfying $1\leq i<j\leq n$ and
$\tau^{-1}\left(  i\right)  >\tau^{-1}\left(  j\right)  $. In other words,
$\left(  \tau\left(  v\right)  ,\tau\left(  u\right)  \right)  $ is an
inversion of $\tau^{-1}$ (by the definition of an \textquotedblleft inversion
of $\tau^{-1}$\textquotedblright). In other words, $\left(  \tau\left(
v\right)  ,\tau\left(  u\right)  \right)  \in\operatorname*{Inv}\left(
\tau^{-1}\right)  $ (since $\operatorname*{Inv}\left(  \tau^{-1}\right)  $ is
the set of all inversions of $\tau^{-1}$).

On the other hand, $\left(  u,v\right)  =c\in\left(  \operatorname*{Inv}%
\tau\right)  \setminus\left(  \operatorname*{Inv}\left(  \sigma\circ
\tau\right)  \right)  $, so that $\left(  u,v\right)  \notin%
\operatorname*{Inv}\left(  \sigma\circ\tau\right)  $. From this, it is easy to
obtain that $\left(  \sigma\circ\tau\right)  \left(  u\right)  \leq\left(
\sigma\circ\tau\right)  \left(  v\right)  $\ \ \ \ \footnote{\textit{Proof.}
Assume the contrary. Thus, we don't have $\left(  \sigma\circ\tau\right)
\left(  u\right)  \leq\left(  \sigma\circ\tau\right)  \left(  v\right)  $.
Hence, we have $\left(  \sigma\circ\tau\right)  \left(  u\right)  >\left(
\sigma\circ\tau\right)  \left(  v\right)  $. Thus, $\left(  u,v\right)  $ is a
pair of integers satisfying $1\leq u<v\leq n$ and $\left(  \sigma\circ
\tau\right)  \left(  u\right)  >\left(  \sigma\circ\tau\right)  \left(
v\right)  $. In other words, $\left(  u,v\right)  $ is a pair $\left(
i,j\right)  $ of integers satisfying $1\leq i<j\leq n$ and $\left(
\sigma\circ\tau\right)  \left(  i\right)  >\left(  \sigma\circ\tau\right)
\left(  j\right)  $. In other words, $\left(  u,v\right)  $ is an inversion of
$\sigma\circ\tau$ (by the definition of an \textquotedblleft inversion of
$\sigma\circ\tau$\textquotedblright). In other words, $\left(  u,v\right)
\in\operatorname*{Inv}\left(  \sigma\circ\tau\right)  $ (since
$\operatorname*{Inv}\left(  \sigma\circ\tau\right)  $ is the set of all
inversions of $\sigma\circ\tau$). This contradicts $\left(  u,v\right)
\notin\operatorname*{Inv}\left(  \sigma\circ\tau\right)  $. This contradiction
proves that our assumption was wrong, qed.}. But the map $\sigma\circ\tau$ is
injective\footnote{\textit{Proof.} We have $\sigma\circ\tau\in S_{n}$. Hence,
$\sigma\circ\tau$ is a permutation of $\left\{  1,2,\ldots,n\right\}  $ (since
$S_{n}$ is the set of all permutations of $\left\{  1,2,\ldots,n\right\}  $).
Thus, $\sigma\circ\tau$ is a bijective map. So the map $\sigma\circ\tau$ is
bijective, and therefore injective. Qed.}. But $u<v$ and therefore $u\neq v$.
Hence, $\left(  \sigma\circ\tau\right)  \left(  u\right)  \neq\left(
\sigma\circ\tau\right)  \left(  v\right)  $ (since the map $\sigma\circ\tau$
is injective). Combining this with $\left(  \sigma\circ\tau\right)  \left(
u\right)  \leq\left(  \sigma\circ\tau\right)  \left(  v\right)  $, we obtain
$\left(  \sigma\circ\tau\right)  \left(  u\right)  <\left(  \sigma\circ
\tau\right)  \left(  v\right)  $. Hence, $\sigma\left(  \tau\left(  u\right)
\right)  =\left(  \sigma\circ\tau\right)  \left(  u\right)  <\left(
\sigma\circ\tau\right)  \left(  v\right)  =\sigma\left(  \tau\left(  v\right)
\right)  $. In other words, $\sigma\left(  \tau\left(  v\right)  \right)
>\sigma\left(  \tau\left(  u\right)  \right)  $.

Now, we know that $\left(  \tau\left(  v\right)  ,\tau\left(  u\right)
\right)  $ is a pair of integers satisfying $1\leq\tau\left(  v\right)
<\tau\left(  u\right)  \leq n$ and $\sigma\left(  \tau\left(  v\right)
\right)  >\sigma\left(  \tau\left(  u\right)  \right)  $. In other words,
$\left(  \tau\left(  v\right)  ,\tau\left(  u\right)  \right)  $ is a pair
$\left(  i,j\right)  $ of integers satisfying $1\leq i<j\leq n$ and
$\sigma\left(  i\right)  >\sigma\left(  j\right)  $. In other words, $\left(
\tau\left(  v\right)  ,\tau\left(  u\right)  \right)  $ is an inversion of
$\sigma$ (by the definition of an \textquotedblleft inversion of $\sigma
$\textquotedblright). In other words, $\left(  \tau\left(  v\right)
,\tau\left(  u\right)  \right)  \in\operatorname*{Inv}\sigma$ (since
$\operatorname*{Inv}\sigma$ is the set of all inversions of $\sigma$).
Combining this with $\left(  \tau\left(  v\right)  ,\tau\left(  u\right)
\right)  \in\operatorname*{Inv}\left(  \tau^{-1}\right)  $, we obtain $\left(
\tau\left(  v\right)  ,\tau\left(  u\right)  \right)  \in\left(
\operatorname*{Inv}\sigma\right)  \cap\left(  \operatorname*{Inv}\left(
\tau^{-1}\right)  \right)  =\varnothing$. Thus, the set $\varnothing$ has at
least one element (namely, the element $\left(  \tau\left(  v\right)
,\tau\left(  u\right)  \right)  $). This contradicts the fact that the set
$\varnothing$ has no elements.

Now, forget that we fixed $c$. We thus have derived a contradiction for each
$c\in\left(  \operatorname*{Inv}\tau\right)  \setminus\left(
\operatorname*{Inv}\left(  \sigma\circ\tau\right)  \right)  $. Hence, there
exists no $c\in\left(  \operatorname*{Inv}\tau\right)  \setminus\left(
\operatorname*{Inv}\left(  \sigma\circ\tau\right)  \right)  $. In other words,
$\left(  \operatorname*{Inv}\tau\right)  \setminus\left(  \operatorname*{Inv}%
\left(  \sigma\circ\tau\right)  \right)  $ is the empty set. In other words,
$\operatorname*{Inv}\tau\subseteq\operatorname*{Inv}\left(  \sigma\circ
\tau\right)  $.

Now, forget our assumption that $\left(  \operatorname*{Inv}\sigma\right)
\cap\left(  \operatorname*{Inv}\left(  \tau^{-1}\right)  \right)
=\varnothing$. We thus have proven that if $\left(  \operatorname*{Inv}%
\sigma\right)  \cap\left(  \operatorname*{Inv}\left(  \tau^{-1}\right)
\right)  =\varnothing$, then $\operatorname*{Inv}\tau\subseteq
\operatorname*{Inv}\left(  \sigma\circ\tau\right)  $. In other words, we have
proven the implication (\ref{sol.perm.Inv.sub.e.2}).]

Combining the two implications (\ref{sol.perm.Inv.sub.e.1}) and
(\ref{sol.perm.Inv.sub.e.2}), we obtain the logical equivalence%
\[
\left(  \operatorname*{Inv}\tau\subseteq\operatorname*{Inv}\left(  \sigma
\circ\tau\right)  \right)  \ \Longleftrightarrow\ \left(  \left(
\operatorname*{Inv}\sigma\right)  \cap\left(  \operatorname*{Inv}\left(
\tau^{-1}\right)  \right)  =\varnothing\right)  .
\]
On the other hand, Exercise \ref{exe.perm.Inv.sub} \textbf{(a)} yields that we
have the logical equivalence%
\[
\left(  \ell\left(  \sigma\circ\tau\right)  =\ell\left(  \sigma\right)
+\ell\left(  \tau\right)  \right)  \ \Longleftrightarrow\ \left(
\operatorname*{Inv}\tau\subseteq\operatorname*{Inv}\left(  \sigma\circ
\tau\right)  \right)  .
\]
Hence, we have the following chain of logical equivalences:%
\begin{align*}
\left(  \ell\left(  \sigma\circ\tau\right)  =\ell\left(  \sigma\right)
+\ell\left(  \tau\right)  \right)  \  &  \Longleftrightarrow\ \left(
\operatorname*{Inv}\tau\subseteq\operatorname*{Inv}\left(  \sigma\circ
\tau\right)  \right) \\
\  &  \Longleftrightarrow\ \left(  \left(  \operatorname*{Inv}\sigma\right)
\cap\left(  \operatorname*{Inv}\left(  \tau^{-1}\right)  \right)
=\varnothing\right)  .
\end{align*}

In other words, $\ell\left(  \sigma\circ\tau\right)  =\ell\left(
\sigma\right)  +\ell\left(  \tau\right)  $ holds if and only if $\left(
\operatorname*{Inv}\sigma\right)  \cap\left(  \operatorname*{Inv}\left(
\tau^{-1}\right)  \right)  =\varnothing$. This solves Exercise
\ref{exe.perm.Inv.sub} \textbf{(e)}.
\end{proof}

\subsection{Solution to Exercise \ref{exe.prod(ai+bi)}}

\begin{vershort}
\begin{proof}
[Solution to Exercise \ref{exe.prod(ai+bi)}.]\textbf{(a)} We shall solve
Exercise \ref{exe.prod(ai+bi)} \textbf{(a)} by induction over $n$:

\textit{Induction base:} Exercise \ref{exe.prod(ai+bi)} \textbf{(a)} holds in
the case when $n=0$\ \ \ \ \footnote{\textit{Proof.} Assume that $n=0$. We
must show that Exercise \ref{exe.prod(ai+bi)} \textbf{(a)} holds in this case.
\par
We have $n=0$ and thus $\left[  n\right]  =\varnothing$. Hence, there is only
one subset of $\left[  n\right]  $ (namely, $\varnothing$). Therefore,%
\begin{align*}
\sum_{I\subseteq\left[  n\right]  }\left(  \prod_{i\in I}a_{i}\right)  \left(
\prod_{i\in\left[  n\right]  \setminus I}b_{i}\right)   &
=\underbrace{\left(  \prod_{i\in\varnothing}a_{i}\right)  }%
_{\substack{=\left(  \text{empty product}\right)  \\=1}}\underbrace{\left(
\prod_{i\in\left[  n\right]  \setminus\varnothing}b_{i}\right)  }%
_{\substack{=\prod_{i\in\varnothing}b_{i}\\\text{(since }\left[  n\right]
\setminus\varnothing=\left[  n\right]  =\varnothing\text{)}}}=\prod
_{i\in\varnothing}b_{i} =\left(  \text{empty product}\right)  =1.
\end{align*}
Comparing this with%
\begin{align*}
\prod_{i=1}^{n}\left(  a_{i}+b_{i}\right)   &  =\prod_{i=1}^{0}\left(
a_{i}+b_{i}\right)  \ \ \ \ \ \ \ \ \ \ \left(  \text{since }n=0\right) \\
&  =\left(  \text{empty product}\right)  =1,
\end{align*}
we obtain $\prod_{i=1}^{n}\left(  a_{i}+b_{i}\right)  =\sum_{I\subseteq\left[
n\right]  }\left(  \prod_{i\in I}a_{i}\right)  \left(  \prod_{i\in\left[
n\right]  \setminus I}b_{i}\right)  $. Thus, Exercise \ref{exe.prod(ai+bi)}
\textbf{(a)} holds in the case when $n=0$.}. This completes the induction base.

\textit{Induction step:} Let $k$ be a positive integer. Assume that Exercise
\ref{exe.prod(ai+bi)} \textbf{(a)} holds in the case when $n=k-1$. We must
show that Exercise \ref{exe.prod(ai+bi)} \textbf{(a)} holds in the case when
$n=k$.

Let $\mathbb{K}$ be a commutative ring. Let $a_{1},a_{2},\ldots,a_{k}$ be $k$
elements of $\mathbb{K}$. Let $b_{1},b_{2},\ldots,b_{k}$ be $k$ further
elements of $\mathbb{K}$. We have assumed that Exercise \ref{exe.prod(ai+bi)}
\textbf{(a)} holds in the case when $n=k-1$. Hence, we can apply Exercise
\ref{exe.prod(ai+bi)} \textbf{(a)} to $k-1$, $\left(  a_{1},a_{2}%
,\ldots,a_{k-1}\right)  $ and $\left(  b_{1},b_{2},\ldots,b_{k-1}\right)  $
instead of $n$, $\left(  a_{1},a_{2},\ldots,a_{n}\right)  $ and $\left(
b_{1},b_{2},\ldots,b_{n}\right)  $. We thus obtain%
\begin{equation}
\prod_{i=1}^{k-1}\left(  a_{i}+b_{i}\right)  =\sum_{I\subseteq\left[
k-1\right]  }\left(  \prod_{i\in I}a_{i}\right)  \left(  \prod_{i\in\left[
k-1\right]  \setminus I}b_{i}\right)  . \label{sol.prod(ai+bi).short.a.iass}%
\end{equation}

Now, $k$ is a positive integer. Thus, we have $k\in\left[  k\right]  $ and
$\left[  k\right]  =\left[  k-1\right]  \cup\left\{  k\right\}  $ and $\left[
k-1\right]  =\left[  k\right]  \setminus\left\{  k\right\}  $.

Let us introduce a notation: For any set $S$, we let $\mathcal{P}\left(
S\right)  $ denote the powerset of $S$ (that is, the set of all subsets of
$S$). Now, we observe the following fact:

\begin{statement}
\textit{Fact 1:} Let $S$ be any set. Let $s\in S$. Then, the map%
\begin{align}
\mathcal{P}\left(  S\setminus\left\{  s\right\}  \right)   &  \rightarrow
\mathcal{P}\left(  S\right)  \setminus\mathcal{P}\left(  S\setminus\left\{
s\right\}  \right)  ,\nonumber\\
U  &  \mapsto U\cup\left\{  s\right\}
\label{sol.prod(ai+bi).short.a.fact1map}%
\end{align}
is well-defined and a bijection.
\end{statement}

Fact 1 is easy to prove\footnote{The idea behind it is that the subsets of $S$
which are \textbf{not} subsets of $S\setminus\left\{  s\right\}  $ are the
subsets of $S$ that contain $s$, and each such subset can be written in the
form $U\cup\left\{  s\right\}  $ for some $U\subseteq S\setminus\left\{
s\right\}  $.
\par
If you want to prove Fact 1 formally, you need to prove two statements:
\par
\begin{enumerate}
\item The map (\ref{sol.prod(ai+bi).short.a.fact1map}) is well-defined (i.e.,
we have $U\cup\left\{  s\right\}  \in\mathcal{P}\left(  S\right)
\setminus\mathcal{P}\left(  S\setminus\left\{  s\right\}  \right)  $ for each
$U\in\mathcal{P}\left(  S\right)  $).
\par
\item This map is bijective.
\end{enumerate}
\par
Proving the first statement is straightforward. The best way to prove the
second statement is to show that the map
(\ref{sol.prod(ai+bi).short.a.fact1map}) has an inverse -- namely, the map
$\mathcal{P}\left(  S\right)  \setminus\mathcal{P}\left(  S\setminus\left\{
s\right\}  \right)  \rightarrow\mathcal{P}\left(  S\setminus\left\{
s\right\}  \right)  ,\ V\mapsto V\setminus\left\{  s\right\}  $. Of course,
you would also have to show that this latter map is well-defined, too.}. We
can apply Fact 1 to $S=\left[  k\right]  $ and $s=k$; we thus conclude that
the map%
\begin{align*}
\mathcal{P}\left(  \left[  k\right]  \setminus\left\{  k\right\}  \right)   &
\rightarrow\mathcal{P}\left(  \left[  k\right]  \right)  \setminus
\mathcal{P}\left(  \left[  k\right]  \setminus\left\{  k\right\}  \right)  ,\\
U  &  \mapsto U\cup\left\{  k\right\}
\end{align*}
is well-defined and a bijection. Since $\left[  k\right]  \setminus\left\{
k\right\}  =\left[  k-1\right]  $, this rewrites as follows: The map%
\begin{align*}
\mathcal{P}\left(  \left[  k-1\right]  \right)   &  \rightarrow\mathcal{P}%
\left(  \left[  k\right]  \right)  \setminus\mathcal{P}\left(  \left[
k-1\right]  \right)  ,\\
U  &  \mapsto U\cup\left\{  k\right\}
\end{align*}
is well-defined and a bijection.

But $\left[  k-1\right]  \subseteq\left[  k\right]  $. Hence, every subset of
$\left[  k-1\right]  $ is a subset of $\left[  k\right]  $.

Let us make another helpful observation:

\begin{statement}
\textit{Fact 2:} Let $I$ be a subset of $\left[  k-1\right]  $. Then,%
\begin{equation}
\prod_{i\in\left[  k\right]  \setminus I}b_{i}=b_{k}\prod_{i\in\left[
k-1\right]  \setminus I}b_{i} \label{sol.prod(ai+bi).short.a.fact2.1}%
\end{equation}
and%
\begin{equation}
\prod_{i\in I\cup\left\{  k\right\}  }a_{i}=a_{k}\prod_{i\in I}a_{i}
\label{sol.prod(ai+bi).short.a.fact2.2}%
\end{equation}
and%
\begin{equation}
\prod_{i\in\left[  k\right]  \setminus\left(  I\cup\left\{  k\right\}
\right)  }b_{i}=\prod_{i\in\left[  k-1\right]  \setminus I}b_{i}.
\label{sol.prod(ai+bi).short.a.fact2.3}%
\end{equation}

\end{statement}

[\textit{Proof of Fact 2:} We have%
\[
\left(  \left[  k\right]  \setminus I\right)  \setminus\left\{  k\right\}
=\left[  k\right]  \setminus\underbrace{\left(  I\cup\left\{  k\right\}
\right)  }_{=\left\{  k\right\}  \cup I}=\left[  k\right]  \setminus\left(
\left\{  k\right\}  \cup I\right)  =\underbrace{\left(  \left[  k\right]
\setminus\left\{  k\right\}  \right)  }_{=\left[  k-1\right]  }\setminus
I=\left[  k-1\right]  \setminus I.
\]

If we had $k\in I$, then we would have $k\in I\subseteq\left[  k-1\right]  $,
which would contradict the fact that $k\notin\left[  k-1\right]  $. Thus, we
cannot have $k\in I$. In other words, we have $k\notin I$. Combining
$k\in\left[  k\right]  $ with $k\notin I$, we obtain $k\in\left[  k\right]
\setminus I$. Hence, we can split off the factor for $i=k$ from the product
$\prod_{i\in\left[  k\right]  \setminus I}b_{i}$. We thus obtain%
\[
\prod_{i\in\left[  k\right]  \setminus I}b_{i}=b_{k}\prod_{i\in\left(  \left[
k\right]  \setminus I\right)  \setminus\left\{  k\right\}  }b_{i}=b_{k}%
\prod_{i\in\left[  k-1\right]  \setminus I}b_{i}%
\]
(since $\left(  \left[  k\right]  \setminus I\right)  \setminus\left\{
k\right\}  =\left[  k-1\right]  \setminus I$). This proves
(\ref{sol.prod(ai+bi).short.a.fact2.1}).

We have $k\notin I$ and thus $\left(  I\cup\left\{  k\right\}  \right)
\setminus\left\{  k\right\}  =I$.

We have $k\in\left\{  k\right\}  \subseteq I\cup\left\{  k\right\}  $. Thus,
we can split off the factor for $i=k$ from the product $\prod_{i\in
I\cup\left\{  k\right\}  }a_{i}$. Thus, we obtain%
\[
\prod_{i\in I\cup\left\{  k\right\}  }a_{i}=a_{k}\prod_{i\in\left(
I\cup\left\{  k\right\}  \right)  \setminus\left\{  k\right\}  }a_{i}%
=a_{k}\prod_{i\in I}a_{i}%
\]
(since $\left(  I\cup\left\{  k\right\}  \right)  \setminus\left\{  k\right\}
=I$). This proves (\ref{sol.prod(ai+bi).short.a.fact2.2}).

Recall that $\left[  k\right]  \setminus\left(  I\cup\left\{  k\right\}
\right)  =\left[  k-1\right]  \setminus I$. Hence, $\prod_{i\in\left[
k\right]  \setminus\left(  I\cup\left\{  k\right\}  \right)  }b_{i}%
=\prod_{i\in\left[  k-1\right]  \setminus I}b_{i}$. This proves
(\ref{sol.prod(ai+bi).short.a.fact2.3}). Thus, the proof of Fact 2 is complete.]

We have%
\begin{align}
&  \underbrace{\sum_{\substack{I\subseteq\left[  k\right]  ;\\I\in
\mathcal{P}\left(  \left[  k-1\right]  \right)  }}}_{\substack{=\sum
_{\substack{I\subseteq\left[  k\right]  ;\\I\subseteq\left[  k-1\right]
}}=\sum_{I\subseteq\left[  k-1\right]  }\\\text{(since every subset of
}\left[  k-1\right]  \\\text{is a subset of }\left[  k\right]  \text{)}%
}}\left(  \prod_{i\in I}a_{i}\right)  \left(  \prod_{i\in\left[  k\right]
\setminus I}b_{i}\right) \nonumber\\
&  =\sum_{I\subseteq\left[  k-1\right]  }\left(  \prod_{i\in I}a_{i}\right)
\underbrace{\left(  \prod_{i\in\left[  k\right]  \setminus I}b_{i}\right)
}_{\substack{=b_{k}\prod_{i\in\left[  k-1\right]  \setminus I}b_{i}\\\text{(by
(\ref{sol.prod(ai+bi).short.a.fact2.1}))}}}=\sum_{I\subseteq\left[
k-1\right]  }\left(  \prod_{i\in I}a_{i}\right)  b_{k}\left(  \prod
_{i\in\left[  k-1\right]  \setminus I}b_{i}\right) \nonumber\\
&  =b_{k}\underbrace{\sum_{I\subseteq\left[  k-1\right]  }\left(  \prod_{i\in
I}a_{i}\right)  \left(  \prod_{i\in\left[  k-1\right]  \setminus I}%
b_{i}\right)  }_{\substack{=\prod_{i=1}^{k-1}\left(  a_{i}+b_{i}\right)
\\\text{(by (\ref{sol.prod(ai+bi).short.a.iass}))}}}=b_{k}\prod_{i=1}%
^{k-1}\left(  a_{i}+b_{i}\right)  \label{sol.prod(ai+bi).short.a.add1}%
\end{align}
and%
\begin{align}
&  \underbrace{\sum_{\substack{I\subseteq\left[  k\right]  ;\\I\notin%
\mathcal{P}\left(  \left[  k-1\right]  \right)  }}}_{\substack{=\sum
_{\substack{I\in\mathcal{P}\left(  \left[  k\right]  \right)  ;\\I\notin%
\mathcal{P}\left(  \left[  k-1\right]  \right)  }}=\sum_{I\in\mathcal{P}%
\left(  \left[  k\right]  \right)  \setminus\mathcal{P}\left(  \left[
k-1\right]  \right)  }}}\left(  \prod_{i\in I}a_{i}\right)  \left(
\prod_{i\in\left[  k\right]  \setminus I}b_{i}\right) \nonumber\\
&  =\sum_{I\in\mathcal{P}\left(  \left[  k\right]  \right)  \setminus
\mathcal{P}\left(  \left[  k-1\right]  \right)  }\left(  \prod_{i\in I}%
a_{i}\right)  \left(  \prod_{i\in\left[  k\right]  \setminus I}b_{i}\right)
\nonumber\\
&  =\sum_{U\in\mathcal{P}\left(  \left[  k-1\right]  \right)  }\left(
\prod_{i\in U\cup\left\{  k\right\}  }a_{i}\right)  \left(  \prod_{i\in\left[
k\right]  \setminus\left(  U\cup\left\{  k\right\}  \right)  }b_{i}\right)
\nonumber\\
&  \ \ \ \ \ \ \ \ \ \ \left(
\begin{array}
[c]{c}%
\text{here, we have substituted }U\cup\left\{  k\right\}  \text{ for }I\text{
in the sum,}\\
\text{since the map }\mathcal{P}\left(  \left[  k-1\right]  \right)
\rightarrow\mathcal{P}\left(  \left[  k\right]  \right)  \setminus
\mathcal{P}\left(  \left[  k-1\right]  \right)  ,\ U\mapsto U\cup\left\{
k\right\} \\
\text{is a bijection}%
\end{array}
\right) \nonumber\\
&  =\sum_{I\in\mathcal{P}\left(  \left[  k-1\right]  \right)  }%
\underbrace{\left(  \prod_{i\in I\cup\left\{  k\right\}  }a_{i}\right)
}_{\substack{=a_{k}\prod_{i\in I}a_{i}\\\text{(by
(\ref{sol.prod(ai+bi).short.a.fact2.2}))}}}\underbrace{\left(  \prod
_{i\in\left[  k\right]  \setminus\left(  I\cup\left\{  k\right\}  \right)
}b_{i}\right)  }_{\substack{=\prod_{i\in\left[  k-1\right]  \setminus I}%
b_{i}\\\text{(by (\ref{sol.prod(ai+bi).short.a.fact2.3}))}}}\nonumber\\
&  \ \ \ \ \ \ \ \ \ \ \left(  \text{here, we have renamed the summation index
}U\text{ as }I\right) \nonumber\\
&  =\sum_{I\subseteq\left[  k-1\right]  }a_{k}\left(  \prod_{i\in I}%
a_{i}\right)  \left(  \prod_{i\in\left[  k-1\right]  \setminus I}b_{i}\right)
\nonumber\\
&  =a_{k}\underbrace{\sum_{I\subseteq\left[  k-1\right]  }\left(  \prod_{i\in
I}a_{i}\right)  \left(  \prod_{i\in\left[  k-1\right]  \setminus I}%
b_{i}\right)  }_{\substack{=\prod_{i=1}^{k-1}\left(  a_{i}+b_{i}\right)
\\\text{(by (\ref{sol.prod(ai+bi).short.a.iass}))}}}=a_{k}\prod_{i=1}%
^{k-1}\left(  a_{i}+b_{i}\right)  . \label{sol.prod(ai+bi).short.a.add2}%
\end{align}

Now, every subset $I$ of $\left[  k\right]  $ satisfies either $I\in
\mathcal{P}\left(  \left[  k-1\right]  \right)  $ or $I\notin\mathcal{P}%
\left(  \left[  k-1\right]  \right)  $ (but not both). Hence,%
\begin{align*}
&  \sum_{I\subseteq\left[  k\right]  }\left(  \prod_{i\in I}a_{i}\right)
\left(  \prod_{i\in\left[  k\right]  \setminus I}b_{i}\right) \\
&  =\underbrace{\sum_{\substack{I\subseteq\left[  k\right]  ;\\I\in
\mathcal{P}\left(  \left[  k-1\right]  \right)  }}\left(  \prod_{i\in I}%
a_{i}\right)  \left(  \prod_{i\in\left[  k\right]  \setminus I}b_{i}\right)
}_{\substack{=b_{k}\prod_{i=1}^{k-1}\left(  a_{i}+b_{i}\right)  \\\text{(by
(\ref{sol.prod(ai+bi).short.a.add1}))}}}+\underbrace{\sum
_{\substack{I\subseteq\left[  k\right]  ;\\I\notin\mathcal{P}\left(  \left[
k-1\right]  \right)  }}\left(  \prod_{i\in I}a_{i}\right)  \left(  \prod
_{i\in\left[  k-1\right]  \setminus I}b_{i}\right)  }_{\substack{=a_{k}%
\prod_{i=1}^{k-1}\left(  a_{i}+b_{i}\right)  \\\text{(by
(\ref{sol.prod(ai+bi).short.a.add2}))}}}\\
&  =b_{k}\prod_{i=1}^{k-1}\left(  a_{i}+b_{i}\right)  +a_{k}\prod_{i=1}%
^{k-1}\left(  a_{i}+b_{i}\right)  =\left(  b_{k}+a_{k}\right)  \prod
_{i=1}^{k-1}\left(  a_{i}+b_{i}\right) \\
&  =\left(  a_{k}+b_{k}\right)  \prod_{i=1}^{k-1}\left(  a_{i}+b_{i}\right)
=\prod_{i=1}^{k}\left(  a_{i}+b_{i}\right)  .
\end{align*}
In other words,%
\begin{equation}
\prod_{i=1}^{k}\left(  a_{i}+b_{i}\right)  =\sum_{I\subseteq\left[  k\right]
}\left(  \prod_{i\in I}a_{i}\right)  \left(  \prod_{i\in\left[  k\right]
\setminus I}b_{i}\right)  . \label{sol.prod(ai+bi).short.a.almostthere}%
\end{equation}

Now, forget that we fixed $\mathbb{K}$, $\left(  a_{1},a_{2},\ldots
,a_{k}\right)  $ and $\left(  b_{1},b_{2},\ldots,b_{k}\right)  $. We thus have
proven (\ref{sol.prod(ai+bi).short.a.almostthere}) for every commutative ring
$\mathbb{K}$, every $k$ elements $a_{1},a_{2},\ldots,a_{k}$ of $\mathbb{K}$,
and every $k$ elements $b_{1},b_{2},\ldots,b_{k}$ of $\mathbb{K}$. In other
words, we have proven that Exercise \ref{exe.prod(ai+bi)} \textbf{(a)} holds
in the case when $n=k$. This completes the induction step. Exercise
\ref{exe.prod(ai+bi)} \textbf{(a)} is thus proven by induction.

\textbf{(b)} Let $a\in\mathbb{K}$, $b\in\mathbb{K}$ and $n\in\mathbb{N}$. We
must prove (\ref{eq.rings.(a+b)**n}).

For every $k\in\left\{  0,1,\ldots,n\right\}  $, we have%
\begin{equation}
\left(  \text{the number of all }I\subseteq\left[  n\right]  \text{ satisfying
}\left\vert I\right\vert =k\right)  =\dbinom{n}{k}
\label{sol.prod(ai+bi).short.b.num}%
\end{equation}
\footnote{\textit{Proof of (\ref{sol.prod(ai+bi).short.b.num}):} Let
$k\in\left\{  0,1,\ldots,n\right\}  $.
\par
Clearly, $\left[  n\right]  $ is an $n$-element set. Thus, Proposition
\ref{prop.binom.subsets} (applied to $n$, $k$ and $\left[  n\right]  $ instead
of $m$, $n$ and $S$) shows that $\dbinom{n}{k}$ is the number of all
$k$-element subsets of $\left[  n\right]  $. In other words,%
\begin{align*}
\dbinom{n}{k}  &  =\left(  \text{the number of all }k\text{-element subsets of
}\left[  n\right]  \right) \\
&  =\left(  \text{the number of all }I\subseteq\left[  n\right]  \text{
satisfying }\left\vert I\right\vert =k\right)  .
\end{align*}
This proves (\ref{sol.prod(ai+bi).short.b.num}).}.

But Exercise \ref{exe.prod(ai+bi)} \textbf{(a)} (applied to $a_{i}=a$ and
$b_{i}=b$) yields%
\[
\prod_{i=1}^{n}\left(  a+b\right)  =\sum_{I\subseteq\left[  n\right]
}\underbrace{\left(  \prod_{i\in I}a\right)  }_{=a^{\left\vert I\right\vert }%
}\underbrace{\left(  \prod_{i\in\left[  n\right]  \setminus I}b\right)
}_{\substack{=b^{\left\vert \left[  n\right]  \setminus I\right\vert
}\\=b^{\left\vert \left[  n\right]  \right\vert -\left\vert I\right\vert
}\\\text{(since }\left\vert \left[  n\right]  \setminus I\right\vert
=\left\vert \left[  n\right]  \right\vert -\left\vert I\right\vert
\\\text{(since }I\subseteq\left[  n\right]  \text{))}}}=\sum_{I\in
\mathcal{P}\left(  \left[  n\right]  \right)  }a^{\left\vert I\right\vert
}b^{\left\vert \left[  n\right]  \right\vert -\left\vert I\right\vert }.
\]
Comparing this with $\prod_{i=1}^{n}\left(  a+b\right)  =\left(  a+b\right)
^{n}$, we obtain%
\begin{align*}
\left(  a+b\right)  ^{n}  &  =\underbrace{\sum_{I\subseteq\left[  n\right]  }%
}_{\substack{=\sum_{k\in\left\{  0,1,\ldots,n\right\}  }\sum
_{\substack{I\subseteq\left[  n\right]  ;\\\left\vert I\right\vert
=k}}\\\text{(since every subset }I\text{ of }\left[  n\right]
\\\text{satisfies }\left\vert I\right\vert \in\left\{  0,1,\ldots,n\right\}
\\\text{(because }\left\vert I\right\vert \leq\left\vert \left[  n\right]
\right\vert =n\text{))}}}a^{\left\vert I\right\vert }b^{\left\vert \left[
n\right]  \right\vert -\left\vert I\right\vert }=\sum_{k\in\left\{
0,1,\ldots,n\right\}  }\sum_{\substack{I\subseteq\left[  n\right]
;\\\left\vert I\right\vert =k}}\underbrace{a^{\left\vert I\right\vert }%
}_{\substack{=a^{k}\\\text{(since }\left\vert I\right\vert =k\text{)}%
}}\underbrace{b^{\left\vert \left[  n\right]  \right\vert -\left\vert
I\right\vert }}_{\substack{=b^{n-k}\\\text{(since }\left\vert \left[
n\right]  \right\vert =n\\\text{and }\left\vert I\right\vert =k\text{)}}}\\
&  =\underbrace{\sum_{k\in\left\{  0,1,\ldots,n\right\}  }}_{=\sum_{k=0}^{n}%
}\underbrace{\sum_{\substack{I\subseteq\left[  n\right]  ;\\\left\vert
I\right\vert =k}}a^{k}b^{n-k}}_{=\left(  \text{the number of all }%
I\subseteq\left[  n\right]  \text{ satisfying }\left\vert I\right\vert
=k\right)  a^{k}b^{n-k}}\\
&  =\sum_{k=0}^{n}\underbrace{\left(  \text{the number of all }I\subseteq
\left[  n\right]  \text{ satisfying }\left\vert I\right\vert =k\right)
}_{\substack{=\dbinom{n}{k}\\\text{(by (\ref{sol.prod(ai+bi).short.b.num}))}%
}}a^{k}b^{n-k}=\sum_{k=0}^{n}\dbinom{n}{k}a^{k}b^{n-k}.
\end{align*}
Thus, (\ref{eq.rings.(a+b)**n}) is proven. This solves Exercise
\ref{exe.prod(ai+bi)} \textbf{(b)}.
\end{proof}
\end{vershort}

\begin{verlong}
Before we start solving Exercise \ref{exe.prod(ai+bi)}, let us introduce a notation:

\begin{definition}
\label{def.sol.prod(ai+bi).powerset}Let $S$ be a set. Then, $\mathcal{P}%
\left(  S\right)  $ will denote the powerset of $S$ (that is, the set of all
subsets of $S$).
\end{definition}

\begin{proposition}
\label{prop.sol.prod(ai+bi).powerset.lem}Let $S$ be a set. Let $s\in S$. Then,
$\mathcal{P}\left(  S\setminus\left\{  s\right\}  \right)  \subseteq
\mathcal{P}\left(  S\right)  $. Furthermore, the map%
\begin{align*}
\mathcal{P}\left(  S\setminus\left\{  s\right\}  \right)   &  \rightarrow
\mathcal{P}\left(  S\right)  \setminus\mathcal{P}\left(  S\setminus\left\{
s\right\}  \right)  ,\\
U  &  \mapsto U\cup\left\{  s\right\}
\end{align*}
is well-defined and a bijection.
\end{proposition}

Proposition \ref{prop.sol.prod(ai+bi).powerset.lem} is an analogue (more
precisely, a \textquotedblleft rougher\textquotedblright\ version) of
Proposition \ref{prop.sol.prop.binom.subsets.Pm.lem} \textbf{(c)}; indeed, it
differs from Proposition \ref{prop.sol.prop.binom.subsets.Pm.lem} \textbf{(c)}
in that it concerns itself with arbitrary subsets rather than $m$-element and
$\left(  m-1\right)  $-element subsets. Unsurprisingly, the proof of
Proposition \ref{prop.sol.prod(ai+bi).powerset.lem} is analogous to that of
Proposition \ref{prop.sol.prop.binom.subsets.Pm.lem}. For the sake of
completeness, let me show this proof:

\begin{proof}
[Proof of Proposition \ref{prop.sol.prod(ai+bi).powerset.lem}.]We know that
$\mathcal{P}\left(  S\setminus\left\{  s\right\}  \right)  $ is the set of all
subsets of $S\setminus\left\{  s\right\}  $ (by the definition of
$\mathcal{P}\left(  S\setminus\left\{  s\right\}  \right)  $). Also,
$\mathcal{P}\left(  S\right)  $ is the set of all subsets of $S$ (by the
definition of $\mathcal{P}\left(  S\right)  $).

For every $U\in\mathcal{P}\left(  S\setminus\left\{  s\right\}  \right)  $, we
have $U\in\mathcal{P}\left(  S\right)  $\ \ \ \ \footnote{\textit{Proof.} Let
$U\in\mathcal{P}\left(  S\setminus\left\{  s\right\}  \right)  $. We must show
that $U\in\mathcal{P}\left(  S\right)  $.
\par
We have $U\in\mathcal{P}\left(  S\setminus\left\{  s\right\}  \right)  $. In
other words, $U$ is a subset of $S\setminus\left\{  s\right\}  $ (since
$\mathcal{P}\left(  S\setminus\left\{  s\right\}  \right)  $ is the set of all
subsets of $S\setminus\left\{  s\right\}  $). Thus, $U\subseteq S\setminus
\left\{  s\right\}  \subseteq S$. Therefore, $U$ is a subset of $S$. In other
words, $U\in\mathcal{P}\left(  S\right)  $ (since $\mathcal{P}\left(
S\right)  $ is the set of all subsets of $S$). Qed.}. In other words,
$\mathcal{P}\left(  S\setminus\left\{  s\right\}  \right)  \subseteq
\mathcal{P}\left(  S\right)  $.

For every $U\in\mathcal{P}\left(  S\setminus\left\{  s\right\}  \right)  $, we
have $U\cup\left\{  s\right\}  \in\mathcal{P}\left(  S\right)  \setminus
\mathcal{P}\left(  S\setminus\left\{  s\right\}  \right)  $%
\ \ \ \ \footnote{\textit{Proof.} Let $U\in\mathcal{P}\left(  S\setminus
\left\{  s\right\}  \right)  $. We must prove that $U\cup\left\{  s\right\}
\in\mathcal{P}\left(  S\right)  \setminus\mathcal{P}\left(  S\setminus\left\{
s\right\}  \right)  $.
\par
Let $V=U\cup\left\{  s\right\}  $.
\par
We have $U\in\mathcal{P}\left(  S\setminus\left\{  s\right\}  \right)  $. In
other words, $U$ is a subset of $S\setminus\left\{  s\right\}  $ (since
$\mathcal{P}\left(  S\setminus\left\{  s\right\}  \right)  $ is the set of all
subsets of $S\setminus\left\{  s\right\}  $). Thus, $U\subseteq S\setminus
\left\{  s\right\}  $. Now,
\[
V=\underbrace{U}_{\subseteq S\setminus\left\{  s\right\}  \subseteq S}%
\cup\underbrace{\left\{  s\right\}  }_{\substack{\subseteq S\\\text{(since
}s\in S\text{)}}}\subseteq S\cup S=S.
\]
In other words, $V$ is a subset of $S$. In other words, $V\in\mathcal{P}%
\left(  S\right)  $ (since $\mathcal{P}\left(  S\right)  $ is the set of all
subsets of $S$).
\par
Now, we shall prove that $V\notin\mathcal{P}\left(  S\setminus\left\{
s\right\}  \right)  $. Indeed, assume the contrary (for the sake of
contradiction). Thus, $V\in\mathcal{P}\left(  S\setminus\left\{  s\right\}
\right)  $. In other words, $V$ is a subset of $S\setminus\left\{  s\right\}
$ (since $\mathcal{P}\left(  S\setminus\left\{  s\right\}  \right)  $ is the
set of all subsets of $S\setminus\left\{  s\right\}  $). In other words,
$V\subseteq S\setminus\left\{  s\right\}  $. Thus, $s\in\left\{  s\right\}
\subseteq U\cup\left\{  s\right\}  =V\subseteq S\setminus\left\{  s\right\}
$, which contradicts $s\notin S\setminus\left\{  s\right\}  $. This
contradiction proves that our assumption was false. Hence, $V\notin%
\mathcal{P}\left(  S\setminus\left\{  s\right\}  \right)  $ is proven.
\par
Combining $V\in\mathcal{P}\left(  S\right)  $ with $V\notin\mathcal{P}\left(
S\setminus\left\{  s\right\}  \right)  $, we obtain $V\in\mathcal{P}\left(
S\right)  \setminus\mathcal{P}\left(  S\setminus\left\{  s\right\}  \right)
$. Thus, $U\cup\left\{  s\right\}  =V\in\mathcal{P}\left(  S\right)
\setminus\mathcal{P}\left(  S\setminus\left\{  s\right\}  \right)  $, qed.}.
Hence, we can define a map%
\[
\alpha:\mathcal{P}\left(  S\setminus\left\{  s\right\}  \right)
\rightarrow\mathcal{P}\left(  S\right)  \setminus\mathcal{P}\left(
S\setminus\left\{  s\right\}  \right)
\]
by
\begin{equation}
\left(  \alpha\left(  U\right)  =U\cup\left\{  s\right\}
\ \ \ \ \ \ \ \ \ \ \text{for every }U\in\mathcal{P}\left(  S\setminus\left\{
s\right\}  \right)  \right)  .
\label{pf.prop.sol.prod(ai+bi).powerset.lem.alpha()=}%
\end{equation}
Consider this map $\alpha$.

For every $V\in\mathcal{P}\left(  S\right)  \setminus\mathcal{P}\left(
S\setminus\left\{  s\right\}  \right)  $, we have $V\setminus\left\{
s\right\}  \in\mathcal{P}\left(  S\setminus\left\{  s\right\}  \right)
$\ \ \ \ \footnote{\textit{Proof.} Let $V\in\mathcal{P}\left(  S\right)
\setminus\mathcal{P}\left(  S\setminus\left\{  s\right\}  \right)  $. We must
prove that $V\setminus\left\{  s\right\}  \in\mathcal{P}\left(  S\setminus
\left\{  s\right\}  \right)  $.
\par
Let $W=V\setminus\left\{  s\right\}  $.
\par
We have $V\in\mathcal{P}\left(  S\right)  \setminus\mathcal{P}\left(
S\setminus\left\{  s\right\}  \right)  \subseteq\mathcal{P}\left(  S\right)
$. In other words, $V$ is a subset of $S$ (since $\mathcal{P}\left(  S\right)
$ is the set of all subsets of $S$). In other words, $V\subseteq S$.
\par
Notice that $W=\underbrace{V}_{\subseteq S}\setminus\left\{  s\right\}
\subseteq S\setminus\left\{  s\right\}  $. In other words, $W$ is a subset of
$S\setminus\left\{  s\right\}  $. In other words, $W\in\mathcal{P}\left(
S\setminus\left\{  s\right\}  \right)  $ (since $\mathcal{P}\left(
S\setminus\left\{  s\right\}  \right)  $ is the set of all subsets of
$S\setminus\left\{  s\right\}  $). Hence, $V\setminus\left\{  s\right\}
=W\in\mathcal{P}\left(  S\setminus\left\{  s\right\}  \right)  $. Qed.}.
Hence, we can define a map%
\[
\beta:\mathcal{P}\left(  S\right)  \setminus\mathcal{P}\left(  S\setminus
\left\{  s\right\}  \right)  \rightarrow\mathcal{P}\left(  S\setminus\left\{
s\right\}  \right)
\]
by%
\[
\left(  \beta\left(  V\right)  =V\setminus\left\{  s\right\}
\ \ \ \ \ \ \ \ \ \ \text{for every }V\in\mathcal{P}\left(  S\right)
\setminus\mathcal{P}\left(  S\setminus\left\{  s\right\}  \right)  \right)  .
\]
Consider this map $\beta$.

We have $\alpha\circ\beta=\operatorname*{id}$\ \ \ \ \footnote{\textit{Proof.}
Let $V\in\mathcal{P}\left(  S\right)  \setminus\mathcal{P}\left(
S\setminus\left\{  s\right\}  \right)  $. We shall show that $\left(
\alpha\circ\beta\right)  \left(  V\right)  =\operatorname*{id}\left(
V\right)  $.
\par
We have $V\notin\mathcal{P}\left(  S\setminus\left\{  s\right\}  \right)  $
(since $V\in\mathcal{P}\left(  S\right)  \setminus\mathcal{P}\left(
S\setminus\left\{  s\right\}  \right)  $).
\par
We have $V\in\mathcal{P}\left(  S\right)  \setminus\mathcal{P}\left(
S\setminus\left\{  s\right\}  \right)  \subseteq\mathcal{P}\left(  S\right)
$. In other words, $V$ is an subset of $S$ (since $\mathcal{P}\left(
S\right)  $ is the set of all subsets of $S$). In other words, $V\subseteq S$.
\par
Now, we shall prove that $s\in V$. Indeed, assume the contrary. Thus, $s\notin
V$. Hence, $V\setminus\left\{  s\right\}  =V$, so that $V=\underbrace{V}%
_{\subseteq S}\setminus\left\{  s\right\}  \subseteq S\setminus\left\{
s\right\}  $. Hence, $V$ is a subset of $S\setminus\left\{  s\right\}  $. In
other words, $V\in\mathcal{P}\left(  S\setminus\left\{  s\right\}  \right)  $
(since $\mathcal{P}\left(  S\setminus\left\{  s\right\}  \right)  $ is the set
of all subsets of $S\setminus\left\{  s\right\}  $). This contradicts
$V\notin\mathcal{P}\left(  S\setminus\left\{  s\right\}  \right)  $. This
contradiction proves that our assumption was wrong.
\par
Hence, $s\in V$ is proven. Thus, $V\cup\left\{  s\right\}  =V$. Now,
\begin{align*}
\left(  \alpha\circ\beta\right)  \left(  V\right)   &  =\alpha\left(
\beta\left(  V\right)  \right)  =\underbrace{\beta\left(  V\right)
}_{\substack{=V\setminus\left\{  s\right\}  \\\text{(by the definition of
}\beta\text{)}}}\cup\left\{  s\right\}  \ \ \ \ \ \ \ \ \ \ \left(  \text{by
the definition of }\alpha\right) \\
&  =\left(  V\setminus\left\{  s\right\}  \right)  \cup\left\{  s\right\}
=V\cup\left\{  s\right\}  =V=\operatorname*{id}\left(  V\right)  .
\end{align*}
\par
Now, forget that we fixed $V$. We thus have shown that $\left(  \alpha
\circ\beta\right)  \left(  V\right)  =\operatorname*{id}\left(  V\right)  $
for every $V\in\mathcal{P}\left(  S\right)  \setminus\mathcal{P}\left(
S\setminus\left\{  s\right\}  \right)  $. In other words, $\alpha\circ
\beta=\operatorname*{id}$. Qed.} and $\beta\circ\alpha=\operatorname*{id}%
$\ \ \ \ \footnote{\textit{Proof.} Let $U\in\mathcal{P}\left(  S\setminus
\left\{  s\right\}  \right)  $. We shall prove that $\left(  \beta\circ
\alpha\right)  \left(  U\right)  =\operatorname*{id}\left(  U\right)  $.
\par
We have $U\in\mathcal{P}\left(  S\setminus\left\{  s\right\}  \right)  $. In
other words, $U$ is an subset of $S\setminus\left\{  s\right\}  $ (since
$\mathcal{P}\left(  S\setminus\left\{  s\right\}  \right)  $ is the set of all
subsets of $S\setminus\left\{  s\right\}  $). Hence, $U\subseteq
S\setminus\left\{  s\right\}  $.
\par
We have $s\in\left\{  s\right\}  $ and thus $s\notin S\setminus\left\{
s\right\}  $. If we had $s\in U$, then we would have $s\in U\subseteq
S\setminus\left\{  s\right\}  $, which would contradict $s\notin
S\setminus\left\{  s\right\}  $. Thus, we cannot have $s\in U$. In other
words, we have $s\notin U$. Hence, $U\setminus\left\{  s\right\}  =U$.
\par
Now,%
\begin{align*}
\left(  \beta\circ\alpha\right)  \left(  U\right)   &  =\beta\left(
\alpha\left(  U\right)  \right)  =\underbrace{\alpha\left(  U\right)
}_{\substack{=U\cup\left\{  s\right\}  \\\text{(by the definition of }%
\alpha\text{)}}}\setminus\left\{  s\right\}  \ \ \ \ \ \ \ \ \ \ \left(
\text{by the definition of }\beta\right) \\
&  =\left(  U\cup\left\{  s\right\}  \right)  \setminus\left\{  s\right\}
=U\setminus\left\{  s\right\}  =U=\operatorname*{id}\left(  U\right)  .
\end{align*}
\par
Now, forget that we fixed $U$. Thus, we have shown that $\left(  \beta
\circ\alpha\right)  \left(  U\right)  =\operatorname*{id}\left(  U\right)  $
for each $U\in\mathcal{P}\left(  S\setminus\left\{  s\right\}  \right)  $. In
other words, $\beta\circ\alpha=\operatorname*{id}$, qed.}. These two
equalities show that the maps $\alpha$ and $\beta$ are mutually inverse. Thus,
the map $\alpha$ is invertible, i.e., is a bijection.

But $\alpha$ is the map
\begin{align*}
\mathcal{P}\left(  S\setminus\left\{  s\right\}  \right)   &  \rightarrow
\mathcal{P}\left(  S\right)  \setminus\mathcal{P}\left(  S\setminus\left\{
s\right\}  \right)  ,\\
U  &  \mapsto U\cup\left\{  s\right\}
\end{align*}
(because $\alpha$ is the map $\mathcal{P}\left(  S\setminus\left\{  s\right\}
\right)  \rightarrow\mathcal{P}\left(  S\right)  \setminus\mathcal{P}\left(
S\setminus\left\{  s\right\}  \right)  $ defined by
(\ref{pf.prop.sol.prod(ai+bi).powerset.lem.alpha()=})).

We know that the map $\alpha$ is well-defined and a bijection. In other words,
the map
\begin{align*}
\mathcal{P}\left(  S\setminus\left\{  s\right\}  \right)   &  \rightarrow
\mathcal{P}\left(  S\right)  \setminus\mathcal{P}\left(  S\setminus\left\{
s\right\}  \right)  ,\\
U  &  \mapsto U\cup\left\{  s\right\}
\end{align*}
is well-defined and a bijection (since $\alpha$ is the map
\begin{align*}
\mathcal{P}\left(  S\setminus\left\{  s\right\}  \right)   &  \rightarrow
\mathcal{P}\left(  S\right)  \setminus\mathcal{P}\left(  S\setminus\left\{
s\right\}  \right)  ,\\
U  &  \mapsto U\cup\left\{  s\right\}
\end{align*}
). This proves Proposition \ref{prop.sol.prod(ai+bi).powerset.lem}.
\end{proof}

\begin{proof}
[Solution to Exercise \ref{exe.prod(ai+bi)}.]\textbf{(a)} We shall prove that%
\begin{equation}
\prod_{i=1}^{k}\left(  a_{i}+b_{i}\right)  =\sum_{I\in\mathcal{P}\left(
\left[  k\right]  \right)  }\left(  \prod_{i\in I}a_{i}\right)  \left(
\prod_{i\in\left[  k\right]  \setminus I}b_{i}\right)
\label{sol.prod(ai+bi).a.goal}%
\end{equation}
for every $k\in\left\{  0,1,\ldots,n\right\}  $.

[\textit{Proof of (\ref{sol.prod(ai+bi).a.goal}):} We shall prove
(\ref{sol.prod(ai+bi).a.goal}) by induction over $k$:

\textit{Induction base:} The equality (\ref{sol.prod(ai+bi).a.goal}) holds for
$k=0$\ \ \ \ \footnote{\textit{Proof.} Assume that $k=0$. We must prove that
the equality (\ref{sol.prod(ai+bi).a.goal}) holds.
\par
We have $k=0$, and thus
\begin{align*}
\left[  k\right]   &  =\left[  0\right]  =\left\{  1,2,\ldots,0\right\}
\ \ \ \ \ \ \ \ \ \ \left(  \text{by the definition of }\left[  0\right]
\right) \\
&  =\varnothing.
\end{align*}
Hence, $\mathcal{P}\left(  \left[  k\right]  \right)  =\mathcal{P}\left(
\varnothing\right)  =\left\{  \varnothing\right\}  $. Thus, the sum
$\sum_{I\in\mathcal{P}\left(  \left[  k\right]  \right)  }\left(  \prod_{i\in
I}a_{i}\right)  \left(  \prod_{i\in\left[  k\right]  \setminus I}b_{i}\right)
$ has only one addend -- namely, the addend for $I=\varnothing$. Therefore,
this sum simplifies as follows:%
\begin{align*}
&  \sum_{I\in\mathcal{P}\left(  \left[  k\right]  \right)  }\left(
\prod_{i\in I}a_{i}\right)  \left(  \prod_{i\in\left[  k\right]  \setminus
I}b_{i}\right) \\
&  =\underbrace{\left(  \prod_{i\in\varnothing}a_{i}\right)  }_{=\left(
\text{empty product}\right)  =1}\left(  \prod_{i\in\left[  k\right]
\setminus\varnothing}b_{i}\right)  =\prod_{i\in\left[  k\right]
\setminus\varnothing}b_{i}=\prod_{i\in\varnothing}b_{i}%
\ \ \ \ \ \ \ \ \ \ \left(  \text{since }\left[  k\right]  \setminus
\varnothing=\left[  k\right]  =\varnothing\right) \\
&  =\left(  \text{empty product}\right)  =1.
\end{align*}
Comparing this with%
\begin{align*}
\prod_{i=1}^{k}\left(  a_{i}+b_{i}\right)   &  =\prod_{i=1}^{0}\left(
a_{i}+b_{i}\right)  \ \ \ \ \ \ \ \ \ \ \left(  \text{since }k=0\right) \\
&  =\left(  \text{empty product}\right)  =1,
\end{align*}
we obtain $\prod_{i=1}^{k}\left(  a_{i}+b_{i}\right)  =\sum_{I\in
\mathcal{P}\left(  \left[  k\right]  \right)  }\left(  \prod_{i\in I}%
a_{i}\right)  \left(  \prod_{i\in\left[  k\right]  \setminus I}b_{i}\right)
$. In other words, the equality (\ref{sol.prod(ai+bi).a.goal}) holds. Qed.}.
This completes the induction base.

\textit{Induction step:} Let $K\in\left\{  0,1,\ldots,n\right\}  $ be
positive. Assume that the equality (\ref{sol.prod(ai+bi).a.goal}) holds for
$k=K-1$. We must prove that the equality (\ref{sol.prod(ai+bi).a.goal}) holds
for $k=K$.

We know that $K$ is a positive element of $\left\{  0,1,\ldots,n\right\}  $.
In other words, $K\in\left\{  1,2,\ldots,n\right\}  $. Hence, $K-1\in\left\{
0,1,\ldots,n-1\right\}  \subseteq\left\{  0,1,\ldots,n\right\}  $.

We also know that the equality (\ref{sol.prod(ai+bi).a.goal}) holds for
$k=K-1$. In other words, we have%
\begin{equation}
\prod_{i=1}^{K-1}\left(  a_{i}+b_{i}\right)  =\sum_{I\in\mathcal{P}\left(
\left[  K-1\right]  \right)  }\left(  \prod_{i\in I}a_{i}\right)  \left(
\prod_{i\in\left[  K-1\right]  \setminus I}b_{i}\right)  .
\label{sol.prod(ai+bi).a.goal.pf.iass}%
\end{equation}

But $\left[  K\right]  =\left\{  1,2,\ldots,K\right\}  $ (by the definition of
$\left[  K\right]  $) and $\left[  K-1\right]  =\left\{  1,2,\ldots
,K-1\right\}  $ (by the definition of $\left[  K-1\right]  $). Thus,%
\begin{align}
\underbrace{\left[  K\right]  }_{=\left\{  1,2,\ldots,K\right\}  }%
\setminus\left\{  K\right\}   &  =\left\{  1,2,\ldots,K\right\}
\setminus\left\{  K\right\} \nonumber\\
&  =\left\{  1,2,\ldots,K-1\right\}  \ \ \ \ \ \ \ \ \ \ \left(  \text{since
}K\text{ is a positive integer}\right) \nonumber\\
&  =\left[  K-1\right]  . \label{sol.prod(ai+bi).a.goal.pf.diff}%
\end{align}
Also, $K$ is a positive integer; hence, $K\in\left\{  1,2,\ldots,K\right\}
=\left[  K\right]  $. Thus, Proposition
\ref{prop.sol.prod(ai+bi).powerset.lem} (applied to $S=\left[  K\right]  $ and
$s=K$) yields the following two facts:

\begin{itemize}
\item We have $\mathcal{P}\left(  \left[  K\right]  \setminus\left\{
K\right\}  \right)  \subseteq\mathcal{P}\left(  \left[  K\right]  \right)  $.

\item The map%
\begin{align*}
\mathcal{P}\left(  \left[  K\right]  \setminus\left\{  K\right\}  \right)   &
\rightarrow\mathcal{P}\left(  \left[  K\right]  \right)  \setminus
\mathcal{P}\left(  \left[  K\right]  \setminus\left\{  K\right\}  \right)  ,\\
U  &  \mapsto U\cup\left\{  K\right\}
\end{align*}
is well-defined and a bijection.
\end{itemize}

\noindent Using (\ref{sol.prod(ai+bi).a.goal.pf.diff}), we can rewrite these
two facts as follows:

\begin{itemize}
\item We have $\mathcal{P}\left(  \left[  K-1\right]  \right)  \subseteq
\mathcal{P}\left(  \left[  K\right]  \right)  $.

\item The map%
\begin{align*}
\mathcal{P}\left(  \left[  K-1\right]  \right)   &  \rightarrow\mathcal{P}%
\left(  \left[  K\right]  \right)  \setminus\mathcal{P}\left(  \left[
K-1\right]  \right)  ,\\
U  &  \mapsto U\cup\left\{  K\right\}
\end{align*}
is well-defined and a bijection.
\end{itemize}

Now, we shall prove a helpful observation:

\begin{statement}
\textit{Observation 1:} Let $I\in\mathcal{P}\left(  \left[  K-1\right]
\right)  $. Then,%
\begin{equation}
\prod_{i\in\left[  K\right]  \setminus I}b_{i}=b_{K}\prod_{i\in\left[
K-1\right]  \setminus I}b_{i} \label{sol.prod(ai+bi).a.goal.pf.term1}%
\end{equation}
and%
\begin{equation}
\prod_{i\in I\cup\left\{  K\right\}  }a_{i}=a_{K}\prod_{i\in I}a_{i}
\label{sol.prod(ai+bi).a.goal.pf.term2}%
\end{equation}
and%
\begin{equation}
\prod_{i\in\left[  K\right]  \setminus\left(  I\cup\left\{  K\right\}
\right)  }b_{i}=\prod_{i\in\left[  K-1\right]  \setminus I}b_{i}.
\label{sol.prod(ai+bi).a.goal.pf.term3}%
\end{equation}

\end{statement}

[\textit{Proof of Observation 1:} Let $I\in\mathcal{P}\left(  \left[
K-1\right]  \right)  $. We must prove (\ref{sol.prod(ai+bi).a.goal.pf.term1}),
(\ref{sol.prod(ai+bi).a.goal.pf.term2}) and
(\ref{sol.prod(ai+bi).a.goal.pf.term3}).

We have $I\in\mathcal{P}\left(  \left[  K-1\right]  \right)  $. In other
words, $I$ is a subset of $\left[  K-1\right]  $ (since $\mathcal{P}\left(
\left[  K-1\right]  \right)  $ is the set of all subsets of $\left[
K-1\right]  $ (by the definition of $\mathcal{P}\left(  \left[  K-1\right]
\right)  $)). In other words, $I\subseteq\left[  K-1\right]  $.

We have%
\[
\left(  \left[  K\right]  \setminus I\right)  \setminus\left\{  K\right\}
=\left[  K\right]  \setminus\underbrace{\left(  I\cup\left\{  K\right\}
\right)  }_{=\left\{  K\right\}  \cup I}=\left[  K\right]  \setminus\left(
\left\{  K\right\}  \cup I\right)  =\underbrace{\left(  \left[  K\right]
\setminus\left\{  K\right\}  \right)  }_{=\left[  K-1\right]  }\setminus
I=\left[  K-1\right]  \setminus I.
\]

If we had $K\in\left[  K-1\right]  $, then we would have $K\in\left[
K-1\right]  =\left\{  1,2,\ldots,K-1\right\}  $ and thus $K\leq K-1$; but this
would contradict $K>K-1$. Hence, we cannot have $K\in\left[  K-1\right]  $. In
other words, we have $K\notin\left[  K-1\right]  $.

If we had $K\in I$, then we would have $K\in I\subseteq\left[  K-1\right]  $,
which would contradict the fact that $K\notin\left[  K-1\right]  $. Thus, we
cannot have $K\in I$. In other words, we have $K\notin I$. Combining
$K\in\left[  K\right]  $ with $K\notin I$, we obtain $K\in\left[  K\right]
\setminus I$. Hence, we can split off the factor for $i=K$ from the product
$\prod_{i\in\left[  K\right]  \setminus I}b_{i}$. We thus obtain%
\[
\prod_{i\in\left[  K\right]  \setminus I}b_{i}=b_{K}\prod_{i\in\left(  \left[
K\right]  \setminus I\right)  \setminus\left\{  K\right\}  }b_{i}=b_{K}%
\prod_{i\in\left[  K-1\right]  \setminus I}b_{i}%
\]
(since $\left(  \left[  K\right]  \setminus I\right)  \setminus\left\{
K\right\}  =\left[  K-1\right]  \setminus I$). This proves
(\ref{sol.prod(ai+bi).a.goal.pf.term1}).

We have $K\notin I$ and thus $\left(  I\cup\left\{  K\right\}  \right)
\setminus\left\{  K\right\}  =I$.

We have $K\in\left\{  K\right\}  \subseteq I\cup\left\{  K\right\}  $. Thus,
we can split off the factor for $i=K$ from the product $\prod_{i\in
I\cup\left\{  K\right\}  }a_{i}$. Thus, we obtain%
\[
\prod_{i\in I\cup\left\{  K\right\}  }a_{i}=a_{K}\prod_{i\in\left(
I\cup\left\{  K\right\}  \right)  \setminus\left\{  K\right\}  }a_{i}%
=a_{K}\prod_{i\in I}a_{i}%
\]
(since $\left(  I\cup\left\{  K\right\}  \right)  \setminus\left\{  K\right\}
=I$). This proves (\ref{sol.prod(ai+bi).a.goal.pf.term2}).

Recall that $\left[  K\right]  \setminus\left(  I\cup\left\{  K\right\}
\right)  =\left[  K-1\right]  \setminus I$. Hence, $\prod_{i\in\left[
K\right]  \setminus\left(  I\cup\left\{  K\right\}  \right)  }b_{i}%
=\prod_{i\in\left[  K-1\right]  \setminus I}b_{i}$. This proves
(\ref{sol.prod(ai+bi).a.goal.pf.term3}). Thus, the proof of Observation 1 is complete.]

Now,%
\begin{align}
&  \underbrace{\sum_{\substack{I\in\mathcal{P}\left(  \left[  K\right]
\right)  ;\\I\in\mathcal{P}\left(  \left[  K-1\right]  \right)  }%
}}_{\substack{=\sum_{I\in\mathcal{P}\left(  \left[  K-1\right]  \right)
}\\\text{(since }\mathcal{P}\left(  \left[  K-1\right]  \right)
\subseteq\mathcal{P}\left(  \left[  K\right]  \right)  \text{)}}}\left(
\prod_{i\in I}a_{i}\right)  \left(  \prod_{i\in\left[  K\right]  \setminus
I}b_{i}\right) \nonumber\\
&  =\sum_{I\in\mathcal{P}\left(  \left[  K-1\right]  \right)  }\left(
\prod_{i\in I}a_{i}\right)  \underbrace{\left(  \prod_{i\in\left[  K\right]
\setminus I}b_{i}\right)  }_{\substack{=b_{K}\prod_{i\in\left[  K-1\right]
\setminus I}b_{i}\\\text{(by (\ref{sol.prod(ai+bi).a.goal.pf.term1}))}}%
}=\sum_{I\in\mathcal{P}\left(  \left[  K-1\right]  \right)  }\left(
\prod_{i\in I}a_{i}\right)  b_{K}\prod_{i\in\left[  K-1\right]  \setminus
I}b_{i}\nonumber\\
&  =b_{K}\underbrace{\sum_{I\in\mathcal{P}\left(  \left[  K-1\right]  \right)
}\left(  \prod_{i\in I}a_{i}\right)  \left(  \prod_{i\in\left[  K-1\right]
\setminus I}b_{i}\right)  }_{\substack{=\prod_{i=1}^{K-1}\left(  a_{i}%
+b_{i}\right)  \\\text{(by (\ref{sol.prod(ai+bi).a.goal.pf.iass}))}}%
}=b_{K}\prod_{i=1}^{K-1}\left(  a_{i}+b_{i}\right)
\label{sol.prod(ai+bi).a.goal.pf.add1}%
\end{align}
and%
\begin{align}
&  \underbrace{\sum_{\substack{I\in\mathcal{P}\left(  \left[  K\right]
\right)  ;\\I\notin\mathcal{P}\left(  \left[  K-1\right]  \right)  }}}%
_{=\sum_{I\in\mathcal{P}\left(  \left[  K\right]  \right)  \setminus
\mathcal{P}\left(  \left[  K-1\right]  \right)  }}\left(  \prod_{i\in I}%
a_{i}\right)  \left(  \prod_{i\in\left[  K\right]  \setminus I}b_{i}\right)
\nonumber\\
&  =\sum_{I\in\mathcal{P}\left(  \left[  K\right]  \right)  \setminus
\mathcal{P}\left(  \left[  K-1\right]  \right)  }\left(  \prod_{i\in I}%
a_{i}\right)  \left(  \prod_{i\in\left[  K\right]  \setminus I}b_{i}\right)
\nonumber\\
&  =\sum_{U\in\mathcal{P}\left(  \left[  K-1\right]  \right)  }\left(
\prod_{i\in U\cup\left\{  K\right\}  }a_{i}\right)  \left(  \prod_{i\in\left[
K\right]  \setminus\left(  U\cup\left\{  K\right\}  \right)  }b_{i}\right)
\nonumber\\
&  \ \ \ \ \ \ \ \ \ \ \left(
\begin{array}
[c]{c}%
\text{here, we have substituted }U\cup\left\{  K\right\}  \text{ for }I\text{
in the sum, since}\\
\text{the map }\mathcal{P}\left(  \left[  K-1\right]  \right)  \rightarrow
\mathcal{P}\left(  \left[  K\right]  \right)  \setminus\mathcal{P}\left(
\left[  K-1\right]  \right)  ,\ U\mapsto U\cup\left\{  K\right\} \\
\text{is a bijection}%
\end{array}
\right) \nonumber\\
&  =\sum_{I\in\mathcal{P}\left(  \left[  K-1\right]  \right)  }%
\underbrace{\left(  \prod_{i\in I\cup\left\{  K\right\}  }a_{i}\right)
}_{\substack{=a_{K}\prod_{i\in I}a_{i}\\\text{(by
(\ref{sol.prod(ai+bi).a.goal.pf.term2}))}}}\underbrace{\left(  \prod
_{i\in\left[  K\right]  \setminus\left(  I\cup\left\{  K\right\}  \right)
}b_{i}\right)  }_{\substack{=\prod_{i\in\left[  K-1\right]  \setminus I}%
b_{i}\\\text{(by (\ref{sol.prod(ai+bi).a.goal.pf.term3}))}}}\nonumber\\
&  \ \ \ \ \ \ \ \ \ \ \left(  \text{here, we have renamed the summation index
}U\text{ as }I\right) \nonumber\\
&  =\sum_{I\in\mathcal{P}\left(  \left[  K-1\right]  \right)  }a_{K}\left(
\prod_{i\in I}a_{i}\right)  \left(  \prod_{i\in\left[  K-1\right]  \setminus
I}b_{i}\right) \nonumber\\
&  =a_{K}\underbrace{\sum_{I\in\mathcal{P}\left(  \left[  K-1\right]  \right)
}\left(  \prod_{i\in I}a_{i}\right)  \left(  \prod_{i\in\left[  K-1\right]
\setminus I}b_{i}\right)  }_{\substack{=\prod_{i=1}^{K-1}\left(  a_{i}%
+b_{i}\right)  \\\text{(by (\ref{sol.prod(ai+bi).a.goal.pf.iass}))}}%
}=a_{K}\prod_{i=1}^{K-1}\left(  a_{i}+b_{i}\right)  .
\label{sol.prod(ai+bi).a.goal.pf.add2}%
\end{align}

Now,%
\begin{align*}
&  \sum_{I\in\mathcal{P}\left(  \left[  K\right]  \right)  }\left(
\prod_{i\in I}a_{i}\right)  \left(  \prod_{i\in\left[  K\right]  \setminus
I}b_{i}\right) \\
&  =\underbrace{\sum_{\substack{I\in\mathcal{P}\left(  \left[  K\right]
\right)  ;\\I\in\mathcal{P}\left(  \left[  K-1\right]  \right)  }}\left(
\prod_{i\in I}a_{i}\right)  \left(  \prod_{i\in\left[  K\right]  \setminus
I}b_{i}\right)  }_{\substack{=b_{K}\prod_{i=1}^{K-1}\left(  a_{i}%
+b_{i}\right)  \\\text{(by (\ref{sol.prod(ai+bi).a.goal.pf.add1}))}%
}}+\underbrace{\sum_{\substack{I\in\mathcal{P}\left(  \left[  K\right]
\right)  ;\\I\notin\mathcal{P}\left(  \left[  K-1\right]  \right)  }}\left(
\prod_{i\in I}a_{i}\right)  \left(  \prod_{i\in\left[  K\right]  \setminus
I}b_{i}\right)  }_{\substack{=a_{K}\prod_{i=1}^{K-1}\left(  a_{i}%
+b_{i}\right)  \\\text{(by (\ref{sol.prod(ai+bi).a.goal.pf.add2}))}}}\\
&  \ \ \ \ \ \ \ \ \ \ \left(
\begin{array}
[c]{c}%
\text{since every }I\in\mathcal{P}\left(  \left[  K\right]  \right)  \text{
satisfies either }I\in\mathcal{P}\left(  \left[  K-1\right]  \right) \\
\text{or }I\notin\mathcal{P}\left(  \left[  K-1\right]  \right)  \text{ (but
not both)}%
\end{array}
\right) \\
&  =b_{K}\prod_{i=1}^{K-1}\left(  a_{i}+b_{i}\right)  +a_{K}\prod_{i=1}%
^{K-1}\left(  a_{i}+b_{i}\right)  =\underbrace{\left(  b_{K}+a_{K}\right)
}_{=a_{K}+b_{K}}\prod_{i=1}^{K-1}\left(  a_{i}+b_{i}\right) \\
&  =\left(  a_{K}+b_{K}\right)  \prod_{i=1}^{K-1}\left(  a_{i}+b_{i}\right)  .
\end{align*}
Comparing this with%
\[
\prod_{i=1}^{K}\left(  a_{i}+b_{i}\right)  =\left(  a_{K}+b_{K}\right)
\prod_{i=1}^{K-1}\left(  a_{i}+b_{i}\right)
\]
(here, we have split off the factor for $i=K$ from the product, since
$K\in\left\{  1,2,\ldots,K\right\}  $), we obtain%
\[
\prod_{i=1}^{K}\left(  a_{i}+b_{i}\right)  =\sum_{I\in\mathcal{P}\left(
\left[  K\right]  \right)  }\left(  \prod_{i\in I}a_{i}\right)  \left(
\prod_{i\in\left[  K\right]  \setminus I}b_{i}\right)  .
\]
In other words, the equality (\ref{sol.prod(ai+bi).a.goal}) holds for $k=K$.
This completes the induction step. Thus, the induction proof of
(\ref{sol.prod(ai+bi).a.goal}) is complete.]

Now, the equality (\ref{sol.prod(ai+bi).a.goal}) (applied to $k=n$) yields%
\[
\prod_{i=1}^{n}\left(  a_{i}+b_{i}\right)  =\underbrace{\sum_{I\in
\mathcal{P}\left(  \left[  n\right]  \right)  }}_{=\sum_{I\subseteq\left[
n\right]  }}\left(  \prod_{i\in I}a_{i}\right)  \left(  \prod_{i\in\left[
n\right]  \setminus I}b_{i}\right)  =\sum_{I\subseteq\left[  n\right]
}\left(  \prod_{i\in I}a_{i}\right)  \left(  \prod_{i\in\left[  n\right]
\setminus I}b_{i}\right)  .
\]
This solves Exercise \ref{exe.prod(ai+bi)} \textbf{(a)}.

\textbf{(b)} Let $a\in\mathbb{K}$, $b\in\mathbb{K}$ and $n\in\mathbb{N}$. We
must prove (\ref{eq.rings.(a+b)**n}).

The set $\mathcal{P}\left(  \left[  n\right]  \right)  $ is the set of all
subsets of $\left[  n\right]  $ (by the definition of $\mathcal{P}\left(
\left[  n\right]  \right)  $). Thus, the elements of $\mathcal{P}\left(
\left[  n\right]  \right)  $ are precisely the subsets of $\left[  n\right]  $.

Clearly, $\left\vert \underbrace{\left[  n\right]  }_{=\left\{  1,2,\ldots
,n\right\}  }\right\vert =\left\vert \left\{  1,2,\ldots,n\right\}
\right\vert =n$. In other words, $\left[  n\right]  $ is an $n$-element set.

Every $I\in\mathcal{P}\left(  \left[  n\right]  \right)  $ satisfies
$\left\vert I\right\vert \in\left\{  0,1,\ldots,n\right\}  $%
\ \ \ \ \footnote{\textit{Proof.} Let $I\in\mathcal{P}\left(  \left[
n\right]  \right)  $. Thus, $I$ is a subset of $\left[  n\right]  $ (since
$\mathcal{P}\left(  \left[  n\right]  \right)  $ is the set of all subsets of
$\left[  n\right]  $). Thus, $I$ is a finite set (since $\left[  n\right]  $
is a finite set). Hence, $\left\vert I\right\vert \in\mathbb{N}$.
\par
But $I$ is a subset of $\left[  n\right]  $. In other words, $I\subseteq
\left[  n\right]  $, so that $\left\vert I\right\vert \leq\left\vert \left[
n\right]  \right\vert =n$. Combining this with $\left\vert I\right\vert
\in\mathbb{N}$, we obtain $\left\vert I\right\vert \in\left\{  0,1,\ldots
,n\right\}  $. Qed.}. Every $k\in\left\{  0,1,\ldots,n\right\}  $ satisfies%
\begin{equation}
\left\vert \left\{  I\in\mathcal{P}\left(  \left[  n\right]  \right)
\ \mid\ \left\vert I\right\vert =k\right\}  \right\vert =\dbinom{n}{k}
\label{sol.prod(ai+bi).b.lem}%
\end{equation}
\footnote{\textit{Proof of (\ref{sol.prod(ai+bi).b.lem}):} Let $k\in\left\{
0,1,\ldots,n\right\}  $. Then, Proposition \ref{prop.binom.subsets} (applied
to $n$, $k$ and $\left[  n\right]  $ instead of $m$, $n$ and $S$) shows that
$\dbinom{n}{k}$ is the number of all $k$-element subsets of $\left[  n\right]
$. In other words,%
\[
\dbinom{n}{k}=\left(  \text{the number of all }k\text{-element subsets of
}\left[  n\right]  \right)  .
\]
Comparing this with%
\begin{align*}
&  \left\vert \left\{  I\in\mathcal{P}\left(  \left[  n\right]  \right)
\ \mid\ \left\vert I\right\vert =k\right\}  \right\vert \\
&  =\left(  \text{the number of all }I\in\mathcal{P}\left(  \left[  n\right]
\right)  \text{ satisfying }\left\vert I\right\vert =k\right) \\
&  =\left(  \text{the number of all subsets }I\text{ of }\left[  n\right]
\text{ satisfying }\left\vert I\right\vert =k\right) \\
&  \ \ \ \ \ \ \ \ \ \ \left(  \text{since the elements of }\mathcal{P}\left(
\left[  n\right]  \right)  \text{ are precisely the subsets of }\left[
n\right]  \right) \\
&  =\left(  \text{the number of all }k\text{-element subsets of }\left[
n\right]  \right) \\
&  \ \ \ \ \ \ \ \ \ \ \left(
\begin{array}
[c]{c}%
\text{since the subsets }I\text{ of }\left[  n\right]  \text{ satisfying
}\left\vert I\right\vert =k\text{ are}\\
\text{precisely the }k\text{-element subsets of }\left[  n\right]
\end{array}
\right)  ,
\end{align*}
we obtain $\left\vert \left\{  I\in\mathcal{P}\left(  \left[  n\right]
\right)  \ \mid\ \left\vert I\right\vert =k\right\}  \right\vert =\dbinom
{n}{k}$. This proves (\ref{sol.prod(ai+bi).b.lem}).}.

Exercise \ref{exe.prod(ai+bi)} \textbf{(a)} (applied to $a_{i}=a$ and
$b_{i}=b$) yields%
\[
\prod_{i=1}^{n}\left(  a+b\right)  =\underbrace{\sum_{I\subseteq\left[
n\right]  }}_{=\sum_{I\in\mathcal{P}\left(  \left[  n\right]  \right)  }%
}\underbrace{\left(  \prod_{i\in I}a\right)  }_{=a^{\left\vert I\right\vert }%
}\underbrace{\left(  \prod_{i\in\left[  n\right]  \setminus I}b\right)
}_{\substack{=b^{\left\vert \left[  n\right]  \setminus I\right\vert
}\\=b^{\left\vert \left[  n\right]  \right\vert -\left\vert I\right\vert
}\\\text{(since }\left\vert \left[  n\right]  \setminus I\right\vert
=\left\vert \left[  n\right]  \right\vert -\left\vert I\right\vert
\\\text{(since }I\subseteq\left[  n\right]  \text{))}}}=\sum_{I\in
\mathcal{P}\left(  \left[  n\right]  \right)  }a^{\left\vert I\right\vert
}b^{\left\vert \left[  n\right]  \right\vert -\left\vert I\right\vert }.
\]
Comparing this with $\prod_{i=1}^{n}\left(  a+b\right)  =\left(  a+b\right)
^{n}$, we obtain%
\begin{align*}
\left(  a+b\right)  ^{n}  &  =\underbrace{\sum_{I\in\mathcal{P}\left(  \left[
n\right]  \right)  }}_{\substack{=\sum_{k\in\left\{  0,1,\ldots,n\right\}
}\sum_{\substack{I\in\mathcal{P}\left(  \left[  n\right]  \right)
;\\\left\vert I\right\vert =k}}\\\text{(since every }I\in\mathcal{P}\left(
\left[  n\right]  \right)  \\\text{satisfies }\left\vert I\right\vert
\in\left\{  0,1,\ldots,n\right\}  \text{)}}}a^{\left\vert I\right\vert
}\underbrace{b^{\left\vert \left[  n\right]  \right\vert -\left\vert
I\right\vert }}_{\substack{=b^{n-\left\vert I\right\vert }\\\text{(since
}\left\vert \left[  n\right]  \right\vert =n\text{)}}}=\sum_{k\in\left\{
0,1,\ldots,n\right\}  }\sum_{\substack{I\in\mathcal{P}\left(  \left[
n\right]  \right)  ;\\\left\vert I\right\vert =k}}\underbrace{a^{\left\vert
I\right\vert }}_{\substack{=a^{k}\\\text{(since }\left\vert I\right\vert
=k\text{)}}}\underbrace{b^{n-\left\vert I\right\vert }}_{\substack{=b^{n-k}%
\\\text{(since }\left\vert I\right\vert =k\text{)}}}\\
&  =\underbrace{\sum_{k\in\left\{  0,1,\ldots,n\right\}  }}_{=\sum_{k=0}^{n}%
}\underbrace{\sum_{\substack{I\in\mathcal{P}\left(  \left[  n\right]  \right)
;\\\left\vert I\right\vert =k}}a^{k}b^{n-k}}_{=\left\vert \left\{
I\in\mathcal{P}\left(  \left[  n\right]  \right)  \ \mid\ \left\vert
I\right\vert =k\right\}  \right\vert a^{k}b^{n-k}}\\
&  =\sum_{k=0}^{n}\underbrace{\left\vert \left\{  I\in\mathcal{P}\left(
\left[  n\right]  \right)  \ \mid\ \left\vert I\right\vert =k\right\}
\right\vert }_{\substack{=\dbinom{n}{k}\\\text{(by
(\ref{sol.prod(ai+bi).b.lem}))}}}a^{k}b^{n-k}=\sum_{k=0}^{n}\dbinom{n}{k}%
a^{k}b^{n-k}.
\end{align*}
Thus, (\ref{eq.rings.(a+b)**n}) is proven. This solves Exercise
\ref{exe.prod(ai+bi)} \textbf{(b)}.
\end{proof}
\end{verlong}

\subsection{Solution to Exercise \ref{exe.multinom2}}

Our solution to Exercise \ref{exe.multinom2} will be somewhat similar to our
solution to Exercise \ref{exe.multichoose} above; in particular, it will rely
on Corollary \ref{cor.prodrule.prod-assMZ} again.

In this section, we shall use the notation $\mathbf{m}\left(  k_{1}%
,k_{2},\ldots,k_{m}\right)  $ as defined in Exercise \ref{exe.multinom2}.

Before we solve Exercise \ref{exe.multinom2}, let us record some really
trivial facts:

\begin{lemma}
\label{lem.sol.multinom2.0-tup-m}Each $0$-tuple $\left(  k_{1},k_{2}%
,\ldots,k_{0}\right)  \in\mathbb{N}^{0}$ satisfies $\mathbf{m}\left(
k_{1},k_{2},\ldots,k_{0}\right)  =1$.
\end{lemma}

(Note that there exists only one $0$-tuple $\left(  k_{1},k_{2},\ldots
,k_{0}\right)  \in\mathbb{N}^{0}$, namely the empty list $\left(  {}\right)
$. Thus, the word \textquotedblleft each\textquotedblright\ in Lemma
\ref{lem.sol.multinom2.0-tup-m} is somewhat misleading.)

\begin{proof}
[Proof of Lemma \ref{lem.sol.multinom2.0-tup-m}.]Let $\left(  k_{1}%
,k_{2},\ldots,k_{0}\right)  \in\mathbb{N}^{0}$ be a $0$-tuple. Then, the
definition of $\mathbf{m}\left(  k_{1},k_{2},\ldots,k_{0}\right)  $ yields%
\begin{align*}
\mathbf{m}\left(  k_{1},k_{2},\ldots,k_{0}\right)   &  =\dfrac{\left(
k_{1}+k_{2}+\cdots+k_{0}\right)  !}{k_{1}!k_{2}!\cdots k_{0}!}\\
&  =\left(  \underbrace{k_{1}+k_{2}+\cdots+k_{0}}_{=\left(  \text{empty
sum}\right)  =0}\right)  !/\left(  \underbrace{k_{1}!k_{2}!\cdots k_{0}%
!}_{=\left(  \text{empty product}\right)  =1}\right) \\
&  =0!/1=0!=1.
\end{align*}
This proves Lemma \ref{lem.sol.multinom2.0-tup-m}.
\end{proof}

\begin{lemma}
\label{lem.sol.multinom2.prod}Let $M$ be a positive integer. Let $\left(
s_{1},s_{2},\ldots,s_{M}\right)  \in\mathbb{N}^{M}$ and $n\in\mathbb{N}$ be
such that $s_{1}+s_{2}+\cdots+s_{M}=n$. Then,%
\[
\dbinom{n}{s_{M}}\mathbf{m}\left(  s_{1},s_{2},\ldots,s_{M-1}\right)
=\mathbf{m}\left(  s_{1},s_{2},\ldots,s_{M}\right)  .
\]

\end{lemma}

\begin{proof}
[Proof of Lemma \ref{lem.sol.multinom2.prod}.]We have%
\[
\left(  s_{1}+s_{2}+\cdots+s_{M-1}\right)  +s_{M}=s_{1}+s_{2}+\cdots+s_{M}=n.
\]
Subtracting $s_{M}$ from both sides of this equality, we obtain $s_{1}%
+s_{2}+\cdots+s_{M-1}=n-s_{M}$.

The definition of $\mathbf{m}\left(  s_{1},s_{2},\ldots,s_{M-1}\right)  $
yields%
\begin{equation}
\mathbf{m}\left(  s_{1},s_{2},\ldots,s_{M-1}\right)  =\dfrac{\left(
s_{1}+s_{2}+\cdots+s_{M-1}\right)  !}{s_{1}!s_{2}!\cdots s_{M-1}!}%
=\dfrac{\left(  n-s_{M}\right)  !}{s_{1}!s_{2}!\cdots s_{M-1}!}
\label{pf.lem.sol.multinom2.prod.1}%
\end{equation}
(since $s_{1}+s_{2}+\cdots+s_{M-1}=n-s_{M}$).

The definition of $\mathbf{m}\left(  s_{1},s_{2},\ldots,s_{M}\right)  $ yields%
\begin{equation}
\mathbf{m}\left(  s_{1},s_{2},\ldots,s_{M}\right)  =\dfrac{\left(  s_{1}%
+s_{2}+\cdots+s_{M}\right)  !}{s_{1}!s_{2}!\cdots s_{M}!}=\dfrac{n!}%
{s_{1}!s_{2}!\cdots s_{M}!} \label{pf.lem.sol.multinom2.prod.2}%
\end{equation}
(since $s_{1}+s_{2}+\cdots+s_{M}=n$).

\begin{vershort}
On the other hand, $n=\underbrace{\left(  s_{1}+s_{2}+\cdots+s_{M-1}\right)
}_{\geq0}+s_{M}\geq s_{M}$. Hence, Proposition \ref{prop.binom.formula}
(applied to $n$ and $s_{M}$ instead of $m$ and $n$) yields $\dbinom{n}{s_{M}%
}=\dfrac{n!}{s_{M}!\left(  n-s_{M}\right)  !}$. Multiplying this equality with
(\ref{pf.lem.sol.multinom2.prod.1}), we obtain%
\begin{align*}
&  \dbinom{n}{s_{M}}\mathbf{m}\left(  s_{1},s_{2},\ldots,s_{M-1}\right) \\
&  =\dfrac{n!}{s_{M}!\left(  n-s_{M}\right)  !}\cdot\dfrac{\left(
n-s_{M}\right)  !}{s_{1}!s_{2}!\cdots s_{M-1}!}=\dfrac{n!}{\left(  s_{1}%
!s_{2}!\cdots s_{M-1}!\right)  s_{M}!}=\dfrac{n!}{s_{1}!s_{2}!\cdots s_{M}!}\\
&  =\mathbf{m}\left(  s_{1},s_{2},\ldots,s_{M}\right)
\ \ \ \ \ \ \ \ \ \ \left(  \text{by (\ref{pf.lem.sol.multinom2.prod.2}%
)}\right)  .
\end{align*}
This proves Lemma \ref{lem.sol.multinom2.prod}. \qedhere

\end{vershort}

\begin{verlong}
On the other hand, $s_{1},s_{2},\ldots,s_{M}$ are elements of $\mathbb{N}$
(since $\left(  s_{1},s_{2},\ldots,s_{M}\right)  \in\mathbb{N}^{M}$). Hence,
$s_{1}+s_{2}+\cdots+s_{M-1}\in\mathbb{N}$ and $s_{M}\in\mathbb{N}$.

We have
\[
n=s_{1}+s_{2}+\cdots+s_{M}=\underbrace{\left(  s_{1}+s_{2}+\cdots
+s_{M-1}\right)  }_{\substack{\geq0\\\text{(since }s_{1}+s_{2}+\cdots
+s_{M-1}\in\mathbb{N}\text{)}}}+s_{M}\geq s_{M}.
\]
Hence, Proposition \ref{prop.binom.formula} (applied to $n$ and $s_{M}$
instead of $m$ and $n$) yields $\dbinom{n}{s_{M}}=\dfrac{n!}{s_{M}!\left(
n-s_{M}\right)  !}$. Multiplying this equality with
(\ref{pf.lem.sol.multinom2.prod.1}), we obtain%
\begin{align*}
&  \dbinom{n}{s_{M}}\mathbf{m}\left(  s_{1},s_{2},\ldots,s_{M-1}\right) \\
&  =\dfrac{n!}{s_{M}!\left(  n-s_{M}\right)  !}\cdot\dfrac{\left(
n-s_{M}\right)  !}{s_{1}!s_{2}!\cdots s_{M-1}!}=\dfrac{n!}{s_{M}!}\cdot
\dfrac{1}{s_{1}!s_{2}!\cdots s_{M-1}!}\\
&  =\dfrac{n!}{\left(  s_{1}!s_{2}!\cdots s_{M-1}!\right)  s_{M}!}\\
&  =\dfrac{n!}{s_{1}!s_{2}!\cdots s_{M}!}\ \ \ \ \ \ \ \ \ \ \left(
\text{since }\left(  s_{1}!s_{2}!\cdots s_{M-1}!\right)  s_{M}!=s_{1}%
!s_{2}!\cdots s_{M}!\right) \\
&  =\mathbf{m}\left(  s_{1},s_{2},\ldots,s_{M}\right)
\ \ \ \ \ \ \ \ \ \ \left(  \text{by (\ref{pf.lem.sol.multinom2.prod.2}%
)}\right)  .
\end{align*}
This proves Lemma \ref{lem.sol.multinom2.prod}.
\end{verlong}
\end{proof}

\begin{proof}
[Solution to Exercise \ref{exe.multinom2}.]We shall solve Exercise
\ref{exe.multinom2} by induction over $m$:

\begin{vershort}
\textit{Induction base:} Exercise \ref{exe.multinom2} holds for $m=0$%
\ \ \ \ \footnote{\textit{Proof.} Assume that $m=0$. We must then show that
Exercise \ref{exe.multinom2} holds.
\par
From $m=0$, we obtain $a_{1}+a_{2}+\cdots+a_{m}=a_{1}+a_{2}+\cdots
+a_{0}=\left(  \text{empty sum}\right)  =0$.
\par
We are in one of the following two cases:
\par
\textit{Case 1:} We have $n=0$.
\par
\textit{Case 2:} We have $n\neq0$.
\par
Let us first consider Case 1. In this case, we have $n=0$. Hence, $0^{n}%
=0^{0}=1$. There exists exactly one $0$-tuple $\left(  k_{1},k_{2}%
,\ldots,k_{0}\right)  \in\mathbb{N}^{0}$ (namely, the empty list $\left(
{}\right)  $), and this $0$-tuple $\left(  k_{1},k_{2},\ldots,k_{0}\right)  $
satisfies $k_{1}+k_{2}+\cdots+k_{0}=\left(  \text{empty sum}\right)  =0=n$.
Hence, there exists exactly one $0$-tuple $\left(  k_{1},k_{2},\ldots
,k_{0}\right)  \in\mathbb{N}^{0}$ satisfying $k_{1}+k_{2}+\cdots+k_{0}=n$. The
sum $\sum_{\substack{\left(  k_{1},k_{2},\ldots,k_{0}\right)  \in
\mathbb{N}^{0};\\k_{1}+k_{2}+\cdots+k_{0}=n}}1$ therefore has exactly $1$
addend; thus, this sum rewrites as $\sum_{\substack{\left(  k_{1},k_{2}%
,\ldots,k_{0}\right)  \in\mathbb{N}^{0};\\k_{1}+k_{2}+\cdots+k_{0}=n}}1=1$.
\par
But recall that $m=0$. Hence,%
\begin{align*}
&  \sum_{\substack{\left(  k_{1},k_{2},\ldots,k_{m}\right)  \in\mathbb{N}%
^{m};\\k_{1}+k_{2}+\cdots+k_{m}=n}}\mathbf{m}\left(  k_{1},k_{2},\ldots
,k_{m}\right)  \prod_{i=1}^{m}a_{i}^{k_{i}}\\
&  =\sum_{\substack{\left(  k_{1},k_{2},\ldots,k_{0}\right)  \in\mathbb{N}%
^{0};\\k_{1}+k_{2}+\cdots+k_{0}=n}}\underbrace{\mathbf{m}\left(  k_{1}%
,k_{2},\ldots,k_{0}\right)  }_{\substack{=1\\\text{(by Lemma
\ref{lem.sol.multinom2.0-tup-m})}}}\underbrace{\prod_{i=1}^{0}a_{i}^{k_{i}}%
}_{=\left(  \text{empty product}\right)  =1}=\sum_{\substack{\left(
k_{1},k_{2},\ldots,k_{0}\right)  \in\mathbb{N}^{0};\\k_{1}+k_{2}+\cdots
+k_{0}=n}}1=1.
\end{align*}
Comparing this with%
\[
\left(  \underbrace{a_{1}+a_{2}+\cdots+a_{m}}_{=\left(  \text{empty
sum}\right)  =0}\right)  ^{n}=0^{n}=1,
\]
we obtain%
\[
\left(  a_{1}+a_{2}+\cdots+a_{m}\right)  ^{n}=\sum_{\substack{\left(
k_{1},k_{2},\ldots,k_{m}\right)  \in\mathbb{N}^{m};\\k_{1}+k_{2}+\cdots
+k_{m}=n}}\mathbf{m}\left(  k_{1},k_{2},\ldots,k_{m}\right)  \prod_{i=1}%
^{m}a_{i}^{k_{i}}.
\]
Hence, Exercise \ref{exe.multinom2} holds. We thus have solved Exercise
\ref{exe.multinom2} in Case 1.
\par
Let us now consider Case 2. In this case, we have $n\neq0$. Hence, $n>0$
(since $n\in\mathbb{N}$). Thus, $0^{n}=0$.
\par
Each $0$-tuple $\left(  k_{1},k_{2},\ldots,k_{0}\right)  \in\mathbb{N}^{0}$
satisfies $k_{1}+k_{2}+\cdots+k_{0}=\left(  \text{empty sum}\right)  =0\neq
n$. In other words, no $\left(  k_{1},k_{2},\ldots,k_{0}\right)  \in
\mathbb{N}^{0}$ satisfies $k_{1}+k_{2}+\cdots+k_{0}=n$. Hence, the sum
$\sum_{\substack{\left(  k_{1},k_{2},\ldots,k_{0}\right)  \in\mathbb{N}%
^{0};\\k_{1}+k_{2}+\cdots+k_{0}=n}}\mathbf{m}\left(  k_{1},k_{2},\ldots
,k_{0}\right)  \prod_{i=1}^{0}a_{i}^{k_{i}}$ is an empty sum. Therefore, this
sum rewrites as follows:%
\[
\sum_{\substack{\left(  k_{1},k_{2},\ldots,k_{0}\right)  \in\mathbb{N}%
^{0};\\k_{1}+k_{2}+\cdots+k_{0}=n}}\mathbf{m}\left(  k_{1},k_{2},\ldots
,k_{0}\right)  \prod_{i=1}^{0}a_{i}^{k_{i}}=\left(  \text{empty sum}\right)
=0.
\]
\par
Now, recall that $m=0$. Hence,%
\[
\sum_{\substack{\left(  k_{1},k_{2},\ldots,k_{m}\right)  \in\mathbb{N}%
^{m};\\k_{1}+k_{2}+\cdots+k_{m}=n}}\mathbf{m}\left(  k_{1},k_{2},\ldots
,k_{m}\right)  \prod_{i=1}^{m}a_{i}^{k_{i}}=\sum_{\substack{\left(
k_{1},k_{2},\ldots,k_{0}\right)  \in\mathbb{N}^{0};\\k_{1}+k_{2}+\cdots
+k_{0}=n}}\mathbf{m}\left(  k_{1},k_{2},\ldots,k_{0}\right)  \prod_{i=1}%
^{0}a_{i}^{k_{i}}=0.
\]
Comparing this with%
\[
\left(  \underbrace{a_{1}+a_{2}+\cdots+a_{m}}_{=\left(  \text{empty
sum}\right)  =0}\right)  ^{n}=0^{n}=0,
\]
we obtain%
\[
\left(  a_{1}+a_{2}+\cdots+a_{m}\right)  ^{n}=\sum_{\substack{\left(
k_{1},k_{2},\ldots,k_{m}\right)  \in\mathbb{N}^{m};\\k_{1}+k_{2}+\cdots
+k_{m}=n}}\mathbf{m}\left(  k_{1},k_{2},\ldots,k_{m}\right)  \prod_{i=1}%
^{m}a_{i}^{k_{i}}.
\]
Hence, Exercise \ref{exe.multinom2} holds. We thus have solved Exercise
\ref{exe.multinom2} in Case 2.
\par
We thus have solved Exercise \ref{exe.multinom2} in each of the two Cases 1
and 2. Thus, Exercise \ref{exe.multinom2} always holds (under the assumption
that $m=0$). Qed.}. This completes the induction base.
\end{vershort}

\begin{verlong}
\textit{Induction base:} Exercise \ref{exe.multinom2} holds for $m=0$%
\ \ \ \ \footnote{\textit{Proof.} Assume that $m=0$. We must then show that
Exercise \ref{exe.multinom2} holds.
\par
From $m=0$, we obtain $a_{1}+a_{2}+\cdots+a_{m}=a_{1}+a_{2}+\cdots
+a_{0}=\left(  \text{empty sum}\right)  =0$.
\par
We are in one of the following two cases:
\par
\textit{Case 1:} We have $n=0$.
\par
\textit{Case 2:} We have $n\neq0$.
\par
Let us first consider Case 1. In this case, we have $n=0$. There exists
exactly one $0$-tuple $\left(  k_{1},k_{2},\ldots,k_{0}\right)  \in
\mathbb{N}^{0}$ (namely, the empty list $\left(  {}\right)  $). Hence, the sum
$\sum_{\left(  k_{1},k_{2},\ldots,k_{0}\right)  \in\mathbb{N}^{0}}1$ has
exactly one addend (namely, the addend corresponding to $\left(  k_{1}%
,k_{2},\ldots,k_{0}\right)  =\left(  {}\right)  $). Hence, this sum rewrites
as follows:%
\[
\sum_{\left(  k_{1},k_{2},\ldots,k_{0}\right)  \in\mathbb{N}^{0}}1=1.
\]
But recall that $m=0$. Hence,%
\begin{align*}
&  \sum_{\substack{\left(  k_{1},k_{2},\ldots,k_{m}\right)  \in\mathbb{N}%
^{m};\\k_{1}+k_{2}+\cdots+k_{m}=n}}\mathbf{m}\left(  k_{1},k_{2},\ldots
,k_{m}\right)  \prod_{i=1}^{m}a_{i}^{k_{i}}\\
&  =\underbrace{\sum_{\substack{\left(  k_{1},k_{2},\ldots,k_{0}\right)
\in\mathbb{N}^{0};\\k_{1}+k_{2}+\cdots+k_{0}=n}}}_{\substack{=\sum_{\left(
k_{1},k_{2},\ldots,k_{0}\right)  \in\mathbb{N}^{0}}\\\text{(since each
}\left(  k_{1},k_{2},\ldots,k_{0}\right)  \in\mathbb{N}^{0}\text{
automatically}\\\text{satisfies }k_{1}+k_{2}+\cdots+k_{0}=n\\\text{(because
for each }\left(  k_{1},k_{2},\ldots,k_{0}\right)  \in\mathbb{N}^{0}%
\text{,}\\\text{we have }k_{1}+k_{2}+\cdots+k_{0}=\left(  \text{empty
sum}\right)  =0=n\text{))}}}\underbrace{\mathbf{m}\left(  k_{1},k_{2}%
,\ldots,k_{0}\right)  }_{\substack{=1\\\text{(by Lemma
\ref{lem.sol.multinom2.0-tup-m})}}}\underbrace{\prod_{i=1}^{0}a_{i}^{k_{i}}%
}_{=\left(  \text{empty product}\right)  =1}\\
&  =\sum_{\left(  k_{1},k_{2},\ldots,k_{0}\right)  \in\mathbb{N}^{0}}1=1.
\end{align*}
Comparing this with%
\begin{align*}
\left(  \underbrace{a_{1}+a_{2}+\cdots+a_{m}}_{=\left(  \text{empty
sum}\right)  =0}\right)  ^{n}  &  =0^{n}=0^{0}\ \ \ \ \ \ \ \ \ \ \left(
\text{since }n=0\right) \\
&  =1,
\end{align*}
we obtain%
\[
\left(  a_{1}+a_{2}+\cdots+a_{m}\right)  ^{n}=\sum_{\substack{\left(
k_{1},k_{2},\ldots,k_{m}\right)  \in\mathbb{N}^{m};\\k_{1}+k_{2}+\cdots
+k_{m}=n}}\mathbf{m}\left(  k_{1},k_{2},\ldots,k_{m}\right)  \prod_{i=1}%
^{m}a_{i}^{k_{i}}.
\]
Hence, Exercise \ref{exe.multinom2} holds. We thus have solved Exercise
\ref{exe.multinom2} in Case 1.
\par
Let us now consider Case 2. In this case, we have $n\neq0$. Hence, $n>0$
(since $n\in\mathbb{N}$).
\par
Each $0$-tuple $\left(  k_{1},k_{2},\ldots,k_{0}\right)  \in\mathbb{N}^{0}$
satisfies $k_{1}+k_{2}+\cdots+k_{0}=\left(  \text{empty sum}\right)  =0\neq
n$. In other words, no $\left(  k_{1},k_{2},\ldots,k_{0}\right)  \in
\mathbb{N}^{0}$ satisfies $k_{1}+k_{2}+\cdots+k_{0}=n$. Hence, the sum
$\sum_{\substack{\left(  k_{1},k_{2},\ldots,k_{0}\right)  \in\mathbb{N}%
^{0};\\k_{1}+k_{2}+\cdots+k_{0}=n}}\mathbf{m}\left(  k_{1},k_{2},\ldots
,k_{0}\right)  \prod_{i=1}^{0}a_{i}^{k_{i}}$ is an empty sum. Therefore, this
sum rewrites as follows:%
\[
\sum_{\substack{\left(  k_{1},k_{2},\ldots,k_{0}\right)  \in\mathbb{N}%
^{0};\\k_{1}+k_{2}+\cdots+k_{0}=n}}\mathbf{m}\left(  k_{1},k_{2},\ldots
,k_{0}\right)  \prod_{i=1}^{0}a_{i}^{k_{i}}=\left(  \text{empty sum}\right)
=0.
\]
\par
Now, recall that $m=0$. Hence,%
\begin{align*}
&  \sum_{\substack{\left(  k_{1},k_{2},\ldots,k_{m}\right)  \in\mathbb{N}%
^{m};\\k_{1}+k_{2}+\cdots+k_{m}=n}}\mathbf{m}\left(  k_{1},k_{2},\ldots
,k_{m}\right)  \prod_{i=1}^{m}a_{i}^{k_{i}}\\
&  =\sum_{\substack{\left(  k_{1},k_{2},\ldots,k_{0}\right)  \in\mathbb{N}%
^{0};\\k_{1}+k_{2}+\cdots+k_{0}=n}}\mathbf{m}\left(  k_{1},k_{2},\ldots
,k_{0}\right)  \prod_{i=1}^{0}a_{i}^{k_{i}}=0.
\end{align*}
Comparing this with%
\[
\left(  \underbrace{a_{1}+a_{2}+\cdots+a_{m}}_{=\left(  \text{empty
sum}\right)  =0}\right)  ^{n}=0^{n}=0\ \ \ \ \ \ \ \ \ \ \left(  \text{since
}n>0\right)  ,
\]
we obtain%
\[
\left(  a_{1}+a_{2}+\cdots+a_{m}\right)  ^{n}=\sum_{\substack{\left(
k_{1},k_{2},\ldots,k_{m}\right)  \in\mathbb{N}^{m};\\k_{1}+k_{2}+\cdots
+k_{m}=n}}\mathbf{m}\left(  k_{1},k_{2},\ldots,k_{m}\right)  \prod_{i=1}%
^{m}a_{i}^{k_{i}}.
\]
Hence, Exercise \ref{exe.multinom2} holds. We thus have solved Exercise
\ref{exe.multinom2} in Case 2.
\par
We thus have solved Exercise \ref{exe.multinom2} in each of the two Cases 1
and 2. Thus, Exercise \ref{exe.multinom2} always holds (under the assumption
that $m=0$). Qed.}. This completes the induction base.
\end{verlong}

\textit{Induction step:} Fix a positive integer $M$. Assume that Exercise
\ref{exe.multinom2} holds for $m=M-1$. We now must show that Exercise
\ref{exe.multinom2} holds for $m=M$.

We have assumed that Exercise \ref{exe.multinom2} holds for $m=M-1$. In other
words, the following fact holds:

\begin{statement}
\textit{Fact 1:} Let $\mathbb{K}$ be a commutative ring. Let $a_{1}%
,a_{2},\ldots,a_{M-1}$ be $M-1$ elements of $\mathbb{K}$. Let $n\in\mathbb{N}%
$. Then,%
\[
\left(  a_{1}+a_{2}+\cdots+a_{M-1}\right)  ^{n}=\sum_{\substack{\left(
k_{1},k_{2},\ldots,k_{M-1}\right)  \in\mathbb{N}^{M-1};\\k_{1}+k_{2}%
+\cdots+k_{M-1}=n}}\mathbf{m}\left(  k_{1},k_{2},\ldots,k_{M-1}\right)
\prod_{i=1}^{M-1}a_{i}^{k_{i}}.
\]

\end{statement}

We must show that Exercise \ref{exe.multinom2} holds for $m=M$. In other
words, we must prove the following fact:

\begin{statement}
\textit{Fact 2:} Let $\mathbb{K}$ be a commutative ring. Let $a_{1}%
,a_{2},\ldots,a_{M}$ be $M$ elements of $\mathbb{K}$. Let $n\in\mathbb{N}$.
Then,%
\[
\left(  a_{1}+a_{2}+\cdots+a_{M}\right)  ^{n}=\sum_{\substack{\left(
k_{1},k_{2},\ldots,k_{M}\right)  \in\mathbb{N}^{M};\\k_{1}+k_{2}+\cdots
+k_{M}=n}}\mathbf{m}\left(  k_{1},k_{2},\ldots,k_{M}\right)  \prod_{i=1}%
^{M}a_{i}^{k_{i}}.
\]

\end{statement}

\begin{vershort}
[\textit{Proof of Fact 2:} We have
\[
a_{1}+a_{2}+\cdots+a_{M}=\left(  a_{1}+a_{2}+\cdots+a_{M-1}\right)
+a_{M}=a_{M}+\left(  a_{1}+a_{2}+\cdots+a_{M-1}\right)  .
\]
Taking both sides of this equality to the $n$-th power, we obtain%
\begin{align}
&  \left(  a_{1}+a_{2}+\cdots+a_{M}\right)  ^{n}\nonumber\\
&  =\left(  a_{M}+\left(  a_{1}+a_{2}+\cdots+a_{M-1}\right)  \right)
^{n}\nonumber\\
&  =\sum_{k=0}^{n}\dbinom{n}{k}a_{M}^{k}\left(  a_{1}+a_{2}+\cdots
+a_{M-1}\right)  ^{n-k}\nonumber\\
&  \ \ \ \ \ \ \ \ \ \ \left(  \text{by (\ref{eq.rings.(a+b)**n}) (applied to
}a=a_{M}\text{ and }b=a_{1}+a_{2}+\cdots+a_{M-1}\text{)}\right) \nonumber\\
&  =\sum_{k=0}^{n}\dbinom{n}{k}\left(  a_{1}+a_{2}+\cdots+a_{M-1}\right)
^{n-k}a_{M}^{k}\nonumber\\
&  =\underbrace{\sum_{r=0}^{n}}_{=\sum_{\substack{r\in\mathbb{N};\\r\leq n}%
}}\dbinom{n}{r}\underbrace{\left(  a_{1}+a_{2}+\cdots+a_{M-1}\right)  ^{n-r}%
}_{\substack{=\sum_{\substack{\left(  k_{1},k_{2},\ldots,k_{M-1}\right)
\in\mathbb{N}^{M-1};\\k_{1}+k_{2}+\cdots+k_{M-1}=n-r}}\mathbf{m}\left(
k_{1},k_{2},\ldots,k_{M-1}\right)  \prod_{i=1}^{M-1}a_{i}^{k_{i}}\\\text{(by
Fact 1, applied to }n-r\text{ instead of }r\text{)}}}a_{M}^{r}\nonumber\\
&  \ \ \ \ \ \ \ \ \ \ \left(  \text{here, we have renamed the summation index
}k\text{ as }r\right) \nonumber\\
&  =\sum_{\substack{r\in\mathbb{N};\\r\leq n}}\dbinom{n}{r}\left(
\sum_{\substack{\left(  k_{1},k_{2},\ldots,k_{M-1}\right)  \in\mathbb{N}%
^{M-1};\\k_{1}+k_{2}+\cdots+k_{M-1}=n-r}}\mathbf{m}\left(  k_{1},k_{2}%
,\ldots,k_{M-1}\right)  \prod_{i=1}^{M-1}a_{i}^{k_{i}}\right)  a_{M}%
^{r}\nonumber\\
&  =\sum_{\substack{r\in\mathbb{N};\\r\leq n}}\sum_{\substack{\left(
k_{1},k_{2},\ldots,k_{M-1}\right)  \in\mathbb{N}^{M-1};\\k_{1}+k_{2}%
+\cdots+k_{M-1}=n-r}}\dbinom{n}{r}\mathbf{m}\left(  k_{1},k_{2},\ldots
,k_{M-1}\right)  \left(  \prod_{i=1}^{M-1}a_{i}^{k_{i}}\right)  a_{M}^{r}.
\label{sol.multinom2.short.3}%
\end{align}

On the other hand, each $\left(  \left(  k_{1},k_{2},\ldots,k_{M-1}\right)
,r\right)  \in\mathbb{N}^{M-1}\times\mathbb{N}$ satisfying $k_{1}+k_{2}%
+\cdots+k_{M-1}=n-r$ must automatically satisfy $r\leq n$%
\ \ \ \ \footnote{\textit{Proof.} Let $\left(  \left(  k_{1},k_{2}%
,\ldots,k_{M-1}\right)  ,r\right)  \in\mathbb{N}^{M-1}\times\mathbb{N}$ be
such that $k_{1}+k_{2}+\cdots+k_{M-1}=n-r$. Then, $n-r=k_{1}+k_{2}%
+\cdots+k_{M-1}\geq0$, so that $r\leq n$.}. Hence, we have the following
equality of summation signs:%
\[
\sum_{\substack{\left(  \left(  k_{1},k_{2},\ldots,k_{M-1}\right)  ,r\right)
\in\mathbb{N}^{M-1}\times\mathbb{N};\\k_{1}+k_{2}+\cdots+k_{M-1}=n-r;\\r\leq
n}}=\sum_{\substack{\left(  \left(  k_{1},k_{2},\ldots,k_{M-1}\right)
,r\right)  \in\mathbb{N}^{M-1}\times\mathbb{N};\\k_{1}+k_{2}+\cdots
+k_{M-1}=n-r}}.
\]
Now, we have the following equality of summation signs:%
\begin{align*}
\sum_{\substack{r\in\mathbb{N};\\r\leq n}}\sum_{\substack{\left(  k_{1}%
,k_{2},\ldots,k_{M-1}\right)  \in\mathbb{N}^{M-1};\\k_{1}+k_{2}+\cdots
+k_{M-1}=n-r}}  &  =\sum_{\left(  k_{1},k_{2},\ldots,k_{M-1}\right)
\in\mathbb{N}^{M-1}}\sum_{\substack{r\in\mathbb{N};\\k_{1}+k_{2}%
+\cdots+k_{M-1}=n-r;\\r\leq n}}\\
&  =\sum_{\substack{\left(  \left(  k_{1},k_{2},\ldots,k_{M-1}\right)
,r\right)  \in\mathbb{N}^{M-1}\times\mathbb{N};\\k_{1}+k_{2}+\cdots
+k_{M-1}=n-r;\\r\leq n}}=\sum_{\substack{\left(  \left(  k_{1},k_{2}%
,\ldots,k_{M-1}\right)  ,r\right)  \in\mathbb{N}^{M-1}\times\mathbb{N}%
;\\k_{1}+k_{2}+\cdots+k_{M-1}=n-r}}.
\end{align*}
Hence, (\ref{sol.multinom2.short.3}) becomes%
\begin{align}
&  \left(  a_{1}+a_{2}+\cdots+a_{M}\right)  ^{n}\nonumber\\
&  =\underbrace{\sum_{\substack{r\in\mathbb{N};\\r\leq n}}\sum
_{\substack{\left(  k_{1},k_{2},\ldots,k_{M-1}\right)  \in\mathbb{N}%
^{M-1};\\k_{1}+k_{2}+\cdots+k_{M-1}=n-r}}}_{=\sum_{\substack{\left(  \left(
k_{1},k_{2},\ldots,k_{M-1}\right)  ,r\right)  \in\mathbb{N}^{M-1}%
\times\mathbb{N};\\k_{1}+k_{2}+\cdots+k_{M-1}=n-r}}}\dbinom{n}{r}%
\mathbf{m}\left(  k_{1},k_{2},\ldots,k_{M-1}\right)  \left(  \prod_{i=1}%
^{M-1}a_{i}^{k_{i}}\right)  a_{M}^{r}\nonumber\\
&  =\sum_{\substack{\left(  \left(  k_{1},k_{2},\ldots,k_{M-1}\right)
,r\right)  \in\mathbb{N}^{M-1}\times\mathbb{N};\\k_{1}+k_{2}+\cdots
+k_{M-1}=n-r}}\dbinom{n}{r}\mathbf{m}\left(  k_{1},k_{2},\ldots,k_{M-1}%
\right)  \left(  \prod_{i=1}^{M-1}a_{i}^{k_{i}}\right)  a_{M}^{r}.
\label{sol.multinom2.short.4}%
\end{align}

But Corollary \ref{cor.prodrule.prod-assMZ} (applied to $Z=\mathbb{N}$) shows
that the map%
\begin{align*}
\mathbb{N}^{M}  &  \rightarrow\mathbb{N}^{M-1}\times\mathbb{N},\\
\left(  s_{1},s_{2},\ldots,s_{M}\right)   &  \mapsto\left(  \left(
s_{1},s_{2},\ldots,s_{M-1}\right)  ,s_{M}\right)
\end{align*}
is a bijection. Hence, we can substitute $\left(  \left(  s_{1},s_{2}%
,\ldots,s_{M-1}\right)  ,s_{M}\right)  $ for $\left(  \left(  k_{1}%
,k_{2},\ldots,k_{M-1}\right)  ,r\right)  $ in the sum on the right hand side
of (\ref{sol.multinom2.short.4}). Thus, we obtain%
\begin{align*}
&  \sum_{\substack{\left(  \left(  k_{1},k_{2},\ldots,k_{M-1}\right)
,r\right)  \in\mathbb{N}^{M-1}\times\mathbb{N};\\k_{1}+k_{2}+\cdots
+k_{M-1}=n-r}}\dbinom{n}{r}\mathbf{m}\left(  k_{1},k_{2},\ldots,k_{M-1}%
\right)  \left(  \prod_{i=1}^{M-1}a_{i}^{k_{i}}\right)  a_{M}^{r}\\
&  =\underbrace{\sum_{\substack{\left(  s_{1},s_{2},\ldots,s_{M}\right)
\in\mathbb{N}^{M};\\s_{1}+s_{2}+\cdots+s_{M-1}=n-s_{M}}}}_{\substack{=\sum
_{\substack{\left(  s_{1},s_{2},\ldots,s_{M}\right)  \in\mathbb{N}^{M}%
;\\s_{1}+s_{2}+\cdots+s_{M}=n}}\\\text{(because for any }\left(  s_{1}%
,s_{2},\ldots,s_{M}\right)  \in\mathbb{N}^{M}\text{,}\\\text{the condition
}\left(  s_{1}+s_{2}+\cdots+s_{M-1}=n-s_{M}\right)  \\\text{is equivalent
to}\\\text{the condition }\left(  s_{1}+s_{2}+\cdots+s_{M}=n\right)  \text{)}%
}}\dbinom{n}{s_{M}}\mathbf{m}\left(  s_{1},s_{2},\ldots,s_{M-1}\right)
\underbrace{\left(  \prod_{i=1}^{M-1}a_{i}^{s_{i}}\right)  a_{M}^{s_{M}}%
}_{\substack{=\prod_{i=1}^{M}a_{i}^{s_{i}}}}\\
&  =\sum_{\substack{\left(  s_{1},s_{2},\ldots,s_{M}\right)  \in\mathbb{N}%
^{M};\\s_{1}+s_{2}+\cdots+s_{M}=n}}\underbrace{\dbinom{n}{s_{M}}%
\mathbf{m}\left(  s_{1},s_{2},\ldots,s_{M-1}\right)  }_{\substack{=\mathbf{m}%
\left(  s_{1},s_{2},\ldots,s_{M}\right)  \\\text{(by Lemma
\ref{lem.sol.multinom2.prod})}}}\prod_{i=1}^{M}a_{i}^{s_{i}}\\
&  =\sum_{\substack{\left(  s_{1},s_{2},\ldots,s_{M}\right)  \in\mathbb{N}%
^{M};\\s_{1}+s_{2}+\cdots+s_{M}=n}}\mathbf{m}\left(  s_{1},s_{2},\ldots
,s_{M}\right)  \prod_{i=1}^{M}a_{i}^{s_{i}}=\sum_{\substack{\left(
k_{1},k_{2},\ldots,k_{M}\right)  \in\mathbb{N}^{M};\\k_{1}+k_{2}+\cdots
+k_{M}=n}}\mathbf{m}\left(  k_{1},k_{2},\ldots,k_{M}\right)  \prod_{i=1}%
^{M}a_{i}^{k_{i}}%
\end{align*}
(here, we have renamed the summation index $\left(  s_{1},s_{2},\ldots
,s_{M}\right)  $ as $\left(  k_{1},k_{2},\ldots,k_{M}\right)  $). Hence,
(\ref{sol.multinom2.short.4}) becomes%
\begin{align*}
&  \left(  a_{1}+a_{2}+\cdots+a_{M}\right)  ^{n}\\
&  =\sum_{\substack{\left(  \left(  k_{1},k_{2},\ldots,k_{M-1}\right)
,r\right)  \in\mathbb{N}^{M-1}\times\mathbb{N};\\k_{1}+k_{2}+\cdots
+k_{M-1}=n-r}}\dbinom{n}{r}\mathbf{m}\left(  k_{1},k_{2},\ldots,k_{M-1}%
\right)  \left(  \prod_{i=1}^{M-1}a_{i}^{k_{i}}\right)  a_{M}^{r}\\
&  =\sum_{\substack{\left(  k_{1},k_{2},\ldots,k_{M}\right)  \in\mathbb{N}%
^{M};\\k_{1}+k_{2}+\cdots+k_{M}=n}}\mathbf{m}\left(  k_{1},k_{2},\ldots
,k_{M}\right)  \prod_{i=1}^{M}a_{i}^{k_{i}}.
\end{align*}
This proves Fact 2.]
\end{vershort}

\begin{verlong}
[\textit{Proof of Fact 2:} We have
\[
a_{1}+a_{2}+\cdots+a_{M}=\left(  a_{1}+a_{2}+\cdots+a_{M-1}\right)
+a_{M}=a_{M}+\left(  a_{1}+a_{2}+\cdots+a_{M-1}\right)  .
\]
Taking both sides of this equality to the $n$-th power, we obtain%
\begin{align*}
\left(  a_{1}+a_{2}+\cdots+a_{M}\right)  ^{n}  &  =\left(  a_{M}+\left(
a_{1}+a_{2}+\cdots+a_{M-1}\right)  \right)  ^{n}\\
&  =\sum_{k=0}^{n}\dbinom{n}{k}a_{M}^{k}\left(  a_{1}+a_{2}+\cdots
+a_{M-1}\right)  ^{n-k}%
\end{align*}
(by (\ref{eq.rings.(a+b)**n}) (applied to $a=a_{M}$ and $b=a_{1}+a_{2}%
+\cdots+a_{M-1}$)). Thus,%
\begin{align}
\left(  a_{1}+a_{2}+\cdots+a_{M}\right)  ^{n}  &  =\sum_{k=0}^{n}\dbinom{n}%
{k}a_{M}^{k}\left(  a_{1}+a_{2}+\cdots+a_{M-1}\right)  ^{n-k}\nonumber\\
&  =\sum_{k=0}^{n}\dbinom{n}{k}\left(  a_{1}+a_{2}+\cdots+a_{M-1}\right)
^{n-k}a_{M}^{k}\nonumber\\
&  =\sum_{r=0}^{n}\dbinom{n}{r}\left(  a_{1}+a_{2}+\cdots+a_{M-1}\right)
^{n-r}a_{M}^{r} \label{sol.multinom2.1}%
\end{align}
(here, we have renamed the summation index $k$ as $r$).

But each $r\in\left\{  0,1,\ldots,n\right\}  $ satisfies%
\begin{align}
&  \left(  a_{1}+a_{2}+\cdots+a_{M-1}\right)  ^{n-r}\nonumber\\
&  =\sum_{\substack{\left(  k_{1},k_{2},\ldots,k_{M-1}\right)  \in
\mathbb{N}^{M-1};\\k_{1}+k_{2}+\cdots+k_{M-1}=n-r}}\mathbf{m}\left(
k_{1},k_{2},\ldots,k_{M-1}\right)  \prod_{i=1}^{M-1}a_{i}^{k_{i}}
\label{sol.multinom2.2}%
\end{align}
\footnote{\textit{Proof of (\ref{sol.multinom2.2}):} Let $r\in\left\{
0,1,\ldots,n\right\}  $. Then, $r\leq n$, so that $n-r\in\mathbb{N}$. Hence,
Fact 1 (applied to $n-r$ instead of $n$) yields%
\[
\left(  a_{1}+a_{2}+\cdots+a_{M-1}\right)  ^{n-r}=\sum_{\substack{\left(
k_{1},k_{2},\ldots,k_{M-1}\right)  \in\mathbb{N}^{M-1};\\k_{1}+k_{2}%
+\cdots+k_{M-1}=n-r}}\mathbf{m}\left(  k_{1},k_{2},\ldots,k_{M-1}\right)
\prod_{i=1}^{M-1}a_{i}^{k_{i}}.
\]
Qed.}. Hence, (\ref{sol.multinom2.1}) becomes%
\begin{align}
&  \left(  a_{1}+a_{2}+\cdots+a_{M}\right)  ^{n}\nonumber\\
&  =\underbrace{\sum_{r=0}^{n}}_{=\sum_{\substack{r\in\mathbb{N};\\r\leq n}%
}}\dbinom{n}{r}\underbrace{\left(  a_{1}+a_{2}+\cdots+a_{M-1}\right)  ^{n-r}%
}_{\substack{=\sum_{\substack{\left(  k_{1},k_{2},\ldots,k_{M-1}\right)
\in\mathbb{N}^{M-1};\\k_{1}+k_{2}+\cdots+k_{M-1}=n-r}}\mathbf{m}\left(
k_{1},k_{2},\ldots,k_{M-1}\right)  \prod_{i=1}^{M-1}a_{i}^{k_{i}}\\\text{(by
(\ref{sol.multinom2.2}))}}}a_{M}^{r}\nonumber\\
&  =\sum_{\substack{r\in\mathbb{N};\\r\leq n}}\dbinom{n}{r}\left(
\sum_{\substack{\left(  k_{1},k_{2},\ldots,k_{M-1}\right)  \in\mathbb{N}%
^{M-1};\\k_{1}+k_{2}+\cdots+k_{M-1}=n-r}}\mathbf{m}\left(  k_{1},k_{2}%
,\ldots,k_{M-1}\right)  \prod_{i=1}^{M-1}a_{i}^{k_{i}}\right)  a_{M}%
^{r}\nonumber\\
&  =\sum_{\substack{r\in\mathbb{N};\\r\leq n}}\sum_{\substack{\left(
k_{1},k_{2},\ldots,k_{M-1}\right)  \in\mathbb{N}^{M-1};\\k_{1}+k_{2}%
+\cdots+k_{M-1}=n-r}}\dbinom{n}{r}\mathbf{m}\left(  k_{1},k_{2},\ldots
,k_{M-1}\right)  \left(  \prod_{i=1}^{M-1}a_{i}^{k_{i}}\right)  a_{M}^{r}.
\label{sol.multinom2.3}%
\end{align}

On the other hand, each $\left(  \left(  k_{1},k_{2},\ldots,k_{M-1}\right)
,r\right)  \in\mathbb{N}^{M-1}\times\mathbb{N}$ satisfying $k_{1}+k_{2}%
+\cdots+k_{M-1}=n-r$ must automatically satisfy $r\leq n$%
\ \ \ \ \footnote{\textit{Proof.} Let $\left(  \left(  k_{1},k_{2}%
,\ldots,k_{M-1}\right)  ,r\right)  \in\mathbb{N}^{M-1}\times\mathbb{N}$ be
such that $k_{1}+k_{2}+\cdots+k_{M-1}=n-r$. Then, $\left(  k_{1},k_{2}%
,\ldots,k_{M-1}\right)  \in\mathbb{N}^{M-1}$ (since $\left(  \left(
k_{1},k_{2},\ldots,k_{M-1}\right)  ,r\right)  \in\mathbb{N}^{M-1}%
\times\mathbb{N}$), so that $k_{1},k_{2},\ldots,k_{M-1}$ are elements of
$\mathbb{N}$. Thus, $k_{1}+k_{2}+\cdots+k_{M-1}\in\mathbb{N}$, so that
$k_{1}+k_{2}+\cdots+k_{M-1}\geq0$. Now, $n-r=k_{1}+k_{2}+\cdots+k_{M-1}\geq0$,
so that $r\leq n$. Qed.}. Hence, we have the following equality of summation
signs:%
\[
\sum_{\substack{\left(  \left(  k_{1},k_{2},\ldots,k_{M-1}\right)  ,r\right)
\in\mathbb{N}^{M-1}\times\mathbb{N};\\k_{1}+k_{2}+\cdots+k_{M-1}=n-r;\\r\leq
n}}=\sum_{\substack{\left(  \left(  k_{1},k_{2},\ldots,k_{M-1}\right)
,r\right)  \in\mathbb{N}^{M-1}\times\mathbb{N};\\k_{1}+k_{2}+\cdots
+k_{M-1}=n-r}}.
\]
Now, we have the following equality of summation signs:%
\begin{align*}
\sum_{\substack{r\in\mathbb{N};\\r\leq n}}\sum_{\substack{\left(  k_{1}%
,k_{2},\ldots,k_{M-1}\right)  \in\mathbb{N}^{M-1};\\k_{1}+k_{2}+\cdots
+k_{M-1}=n-r}}  &  =\sum_{\left(  k_{1},k_{2},\ldots,k_{M-1}\right)
\in\mathbb{N}^{M-1}}\sum_{\substack{r\in\mathbb{N};\\k_{1}+k_{2}%
+\cdots+k_{M-1}=n-r;\\r\leq n}}\\
&  =\sum_{\substack{\left(  \left(  k_{1},k_{2},\ldots,k_{M-1}\right)
,r\right)  \in\mathbb{N}^{M-1}\times\mathbb{N};\\k_{1}+k_{2}+\cdots
+k_{M-1}=n-r;\\r\leq n}}=\sum_{\substack{\left(  \left(  k_{1},k_{2}%
,\ldots,k_{M-1}\right)  ,r\right)  \in\mathbb{N}^{M-1}\times\mathbb{N}%
;\\k_{1}+k_{2}+\cdots+k_{M-1}=n-r}}.
\end{align*}
Hence, (\ref{sol.multinom2.3}) becomes%
\begin{align}
&  \left(  a_{1}+a_{2}+\cdots+a_{M}\right)  ^{n}\nonumber\\
&  =\underbrace{\sum_{\substack{r\in\mathbb{N};\\r\leq n}}\sum
_{\substack{\left(  k_{1},k_{2},\ldots,k_{M-1}\right)  \in\mathbb{N}%
^{M-1};\\k_{1}+k_{2}+\cdots+k_{M-1}=n-r}}}_{=\sum_{\substack{\left(  \left(
k_{1},k_{2},\ldots,k_{M-1}\right)  ,r\right)  \in\mathbb{N}^{M-1}%
\times\mathbb{N};\\k_{1}+k_{2}+\cdots+k_{M-1}=n-r}}}\dbinom{n}{r}%
\mathbf{m}\left(  k_{1},k_{2},\ldots,k_{M-1}\right)  \left(  \prod_{i=1}%
^{M-1}a_{i}^{k_{i}}\right)  a_{M}^{r}\nonumber\\
&  =\sum_{\substack{\left(  \left(  k_{1},k_{2},\ldots,k_{M-1}\right)
,r\right)  \in\mathbb{N}^{M-1}\times\mathbb{N};\\k_{1}+k_{2}+\cdots
+k_{M-1}=n-r}}\dbinom{n}{r}\mathbf{m}\left(  k_{1},k_{2},\ldots,k_{M-1}%
\right)  \left(  \prod_{i=1}^{M-1}a_{i}^{k_{i}}\right)  a_{M}^{r}.
\label{sol.multinom2.4}%
\end{align}

But Corollary \ref{cor.prodrule.prod-assMZ} (applied to $Z=\mathbb{N}$) shows
that the map%
\begin{align*}
\mathbb{N}^{M}  &  \rightarrow\mathbb{N}^{M-1}\times\mathbb{N},\\
\left(  s_{1},s_{2},\ldots,s_{M}\right)   &  \mapsto\left(  \left(
s_{1},s_{2},\ldots,s_{M-1}\right)  ,s_{M}\right)
\end{align*}
is a bijection. Hence, we can substitute $\left(  \left(  s_{1},s_{2}%
,\ldots,s_{M-1}\right)  ,s_{M}\right)  $ for $\left(  \left(  k_{1}%
,k_{2},\ldots,k_{M-1}\right)  ,r\right)  $ in the sum on the right hand side
of (\ref{sol.multinom2.4}). Thus, we obtain%
\begin{align*}
&  \sum_{\substack{\left(  \left(  k_{1},k_{2},\ldots,k_{M-1}\right)
,r\right)  \in\mathbb{N}^{M-1}\times\mathbb{N};\\k_{1}+k_{2}+\cdots
+k_{M-1}=n-r}}\dbinom{n}{r}\mathbf{m}\left(  k_{1},k_{2},\ldots,k_{M-1}%
\right)  \left(  \prod_{i=1}^{M-1}a_{i}^{k_{i}}\right)  a_{M}^{r}\\
&  =\underbrace{\sum_{\substack{\left(  s_{1},s_{2},\ldots,s_{M}\right)
\in\mathbb{N}^{M};\\s_{1}+s_{2}+\cdots+s_{M-1}=n-s_{M}}}}_{\substack{=\sum
_{\substack{\left(  s_{1},s_{2},\ldots,s_{M}\right)  \in\mathbb{N}%
^{M};\\\left(  s_{1}+s_{2}+\cdots+s_{M-1}\right)  +s_{M}=n}}\\\text{(because
for any }\left(  s_{1},s_{2},\ldots,s_{M}\right)  \in\mathbb{N}^{M}%
\text{,}\\\text{the condition }\left(  s_{1}+s_{2}+\cdots+s_{M-1}%
=n-s_{M}\right)  \\\text{is equivalent to}\\\text{the condition }\left(
\left(  s_{1}+s_{2}+\cdots+s_{M-1}\right)  +s_{M}=n\right)  \text{)}}%
}\dbinom{n}{s_{M}}\mathbf{m}\left(  s_{1},s_{2},\ldots,s_{M-1}\right)
\underbrace{\left(  \prod_{i=1}^{M-1}a_{i}^{s_{i}}\right)  a_{M}^{s_{M}}%
}_{\substack{=\prod_{i=1}^{M}a_{i}^{s_{i}}\\\text{(since }\prod_{i=1}^{M}%
a_{i}^{s_{i}}=\left(  \prod_{i=1}^{M-1}a_{i}^{s_{i}}\right)  a_{M}^{s_{M}%
}\\\text{(here, we have split off the}\\\text{factor for }i=M\text{ from the
product))}}}\\
&  =\underbrace{\sum_{\substack{\left(  s_{1},s_{2},\ldots,s_{M}\right)
\in\mathbb{N}^{M};\\\left(  s_{1}+s_{2}+\cdots+s_{M-1}\right)  +s_{M}=n}%
}}_{\substack{=\sum_{\substack{\left(  s_{1},s_{2},\ldots,s_{M}\right)
\in\mathbb{N}^{M};\\s_{1}+s_{2}+\cdots+s_{M}=n}}\\\text{(because for any
}\left(  s_{1},s_{2},\ldots,s_{M}\right)  \in\mathbb{N}^{M}\text{,}\\\text{we
have }\left(  s_{1}+s_{2}+\cdots+s_{M-1}\right)  +s_{M}=s_{1}+s_{2}%
+\cdots+s_{M}\text{)}}}\dbinom{n}{s_{M}}\mathbf{m}\left(  s_{1},s_{2}%
,\ldots,s_{M-1}\right)  \prod_{i=1}^{M}a_{i}^{s_{i}}\\
&  =\sum_{\substack{\left(  s_{1},s_{2},\ldots,s_{M}\right)  \in\mathbb{N}%
^{M};\\s_{1}+s_{2}+\cdots+s_{M}=n}}\underbrace{\dbinom{n}{s_{M}}%
\mathbf{m}\left(  s_{1},s_{2},\ldots,s_{M-1}\right)  }_{\substack{=\mathbf{m}%
\left(  s_{1},s_{2},\ldots,s_{M}\right)  \\\text{(by Lemma
\ref{lem.sol.multinom2.prod})}}}\prod_{i=1}^{M}a_{i}^{s_{i}}\\
&  =\sum_{\substack{\left(  s_{1},s_{2},\ldots,s_{M}\right)  \in\mathbb{N}%
^{M};\\s_{1}+s_{2}+\cdots+s_{M}=n}}\mathbf{m}\left(  s_{1},s_{2},\ldots
,s_{M}\right)  \prod_{i=1}^{M}a_{i}^{s_{i}}\\
&  =\sum_{\substack{\left(  k_{1},k_{2},\ldots,k_{M}\right)  \in\mathbb{N}%
^{M};\\k_{1}+k_{2}+\cdots+k_{M}=n}}\mathbf{m}\left(  k_{1},k_{2},\ldots
,k_{M}\right)  \prod_{i=1}^{M}a_{i}^{k_{i}}%
\end{align*}
(here, we have renamed the summation index $\left(  s_{1},s_{2},\ldots
,s_{M}\right)  $ as $\left(  k_{1},k_{2},\ldots,k_{M}\right)  $). Hence,
(\ref{sol.multinom2.4}) becomes%
\begin{align*}
&  \left(  a_{1}+a_{2}+\cdots+a_{M}\right)  ^{n}\\
&  =\sum_{\substack{\left(  \left(  k_{1},k_{2},\ldots,k_{M-1}\right)
,r\right)  \in\mathbb{N}^{M-1}\times\mathbb{N};\\k_{1}+k_{2}+\cdots
+k_{M-1}=n-r}}\dbinom{n}{r}\mathbf{m}\left(  k_{1},k_{2},\ldots,k_{M-1}%
\right)  \left(  \prod_{i=1}^{M-1}a_{i}^{k_{i}}\right)  a_{M}^{r}\\
&  =\sum_{\substack{\left(  k_{1},k_{2},\ldots,k_{M}\right)  \in\mathbb{N}%
^{M};\\k_{1}+k_{2}+\cdots+k_{M}=n}}\mathbf{m}\left(  k_{1},k_{2},\ldots
,k_{M}\right)  \prod_{i=1}^{M}a_{i}^{k_{i}}.
\end{align*}
This proves Fact 2.]
\end{verlong}

But Fact 2 is precisely the statement of Exercise \ref{exe.multinom2} for
$m=M$. Hence, Exercise \ref{exe.multinom2} holds for $m=M$ (since Fact 2 is
proven). This completes the induction step. Thus, Exercise \ref{exe.multinom2}
is proven by induction.
\end{proof}

\subsection{Solution to Exercise \ref{exe.ps4.3}}

\begin{proof}
[Solution to Exercise \ref{exe.ps4.3}.]We first notice a purely combinatorial
fact: For every $\sigma\in S_{n}$ satisfying $\sigma\neq\operatorname*{id}$,%
\begin{equation}
\text{there exists an }i\in\left\{  1,2,\ldots,n\right\}  \text{ such that
}\sigma\left(  i\right)  >i \label{sol.ps4.3.ineq}%
\end{equation}
\footnote{\textit{Proof of (\ref{sol.ps4.3.ineq}):} Let $\sigma\in S_{n}$ be
such that $\sigma\neq\operatorname*{id}$. We need to prove
(\ref{sol.ps4.3.ineq}).
\par
Assume the contrary. Thus, there exists no $i\in\left\{  1,2,\ldots,n\right\}
$ such that $\sigma\left(  i\right)  > i$. In other words, every $i\in\left\{
1,2,\ldots,n\right\}  $ satisfies $\sigma\left(  i\right)  \leq i$.
\par
We shall now show that every $p\in\left\{  1,2,\ldots,n\right\}  $ satisfies%
\begin{equation}
\sigma\left(  p\right)  =p. \label{sol.ps4.3.ineq.pf.goal}%
\end{equation}
\par
\textit{Proof of (\ref{sol.ps4.3.ineq.pf.goal}):} We will prove
(\ref{sol.ps4.3.ineq.pf.goal}) by strong induction over $p$. Thus, fix some
$P\in\left\{  1,2,\ldots,n\right\}  $. We assume that
(\ref{sol.ps4.3.ineq.pf.goal}) is proven for every $p<P$. We need to show that
(\ref{sol.ps4.3.ineq.pf.goal}) holds for $p=P$.
\par
We have assumed that (\ref{sol.ps4.3.ineq.pf.goal}) is proven for every $p<P$.
In other words,%
\begin{equation}
\sigma\left(  p\right)  =p\ \ \ \ \ \ \ \ \ \ \text{for every }p\in\left\{
1,2,\ldots,n\right\}  \text{ satisfying }p<P.
\label{sol.ps4.3.ineq.pf.goal.pf.hyp}%
\end{equation}
\par
Now, we assume (for the sake of contradiction) that $\sigma\left(  P\right)
\neq P$. Recall that every $i\in\left\{  1,2,\ldots,n\right\}  $ satisfies
$\sigma\left(  i\right)  \leq i$. Applying this to $i=P$, we obtain
$\sigma\left(  P\right)  \leq P$. Combined with $\sigma\left(  P\right)  \neq
P$, this yields $\sigma\left(  P\right)  <P$. Hence,
(\ref{sol.ps4.3.ineq.pf.goal.pf.hyp}) (applied to $p=\sigma\left(  P\right)
$) yields $\sigma\left(  \sigma\left(  P\right)  \right)  =\sigma\left(
P\right)  $.
\par
But the map $\sigma$ is a permutation (since $\sigma\in S_{n}$), thus
injective. Hence, from $\sigma\left(  \sigma\left(  P\right)  \right)
=\sigma\left(  P\right)  $, we obtain $\sigma\left(  P\right)  =P$. This
contradicts $\sigma\left(  P\right)  \neq P$. Hence, we have obtained a
contradiction; thus, our assumption (that $\sigma\left(  P\right)  \neq P$)
must have been wrong. We thus have $\sigma\left(  P\right)  =P$. In other
words, (\ref{sol.ps4.3.ineq.pf.goal}) holds for $p=P$. This completes the
inductive proof of (\ref{sol.ps4.3.ineq.pf.goal}).
\par
Now, (\ref{sol.ps4.3.ineq.pf.goal}) shows that every $p\in\left\{
1,2,\ldots,n\right\}  $ satisfies $\sigma\left(  p\right)
=p=\operatorname*{id}\left(  p\right)  $. In other words, $\sigma
=\operatorname*{id}$. This contradicts $\sigma\neq\operatorname*{id}$. This
contradiction proves that our assumption was wrong. Hence,
(\ref{sol.ps4.3.ineq}) is proven.}. Thus, for every $\sigma\in S_{n}$
satisfying $\sigma\neq\operatorname*{id}$, we have%
\begin{equation}
\prod_{i=1}^{n}a_{i,\sigma\left(  i\right)  }=0 \label{sol.ps4.3.0}%
\end{equation}
\footnote{\textit{Proof of (\ref{sol.ps4.3.0}):} Let $\sigma\in S_{n}$ be such
that $\sigma\neq\operatorname*{id}$. According to (\ref{sol.ps4.3.ineq}), we
know that there exists an $i\in\left\{  1,2,\ldots,n\right\}  $ such that
$\sigma\left(  i\right)  >i$. Let $k$ be such an $i$. Thus, $k$ is an element
of $\left\{  1,2,\ldots,n\right\}  $ satisfying $\sigma\left(  k\right)  >k$.
Hence, $k<\sigma\left(  k\right)  $.
\par
Recall that $a_{i,j}=0$ for every $\left(  i,j\right)  \in\left\{
1,2,\ldots,n\right\}  ^{2}$ satisfying $i<j$. Applying this to $i=k$ and
$j=\sigma\left(  k\right)  $, we obtain $a_{k,\sigma\left(  k\right)  }=0$
(since $k<\sigma\left(  k\right)  $). Hence, one factor of the product
$\prod_{i=1}^{n}a_{i,\sigma\left(  i\right)  }$ is $0$ (namely, the factor
$a_{k,\sigma\left(  k\right)  }$). Therefore, the product $\prod_{i=1}%
^{n}a_{i,\sigma\left(  i\right)  }$ is $0$ (because if one factor of a product
is $0$, then the whole product is $0$).}.

Now, (\ref{eq.det.eq.2}) yields%
\begin{align*}
\det A  &  =\sum_{\sigma\in S_{n}}\left(  -1\right)  ^{\sigma}\prod_{i=1}%
^{n}a_{i,\sigma\left(  i\right)  }\\
&  =\underbrace{\left(  -1\right)  ^{\operatorname*{id}}}_{=1}\prod_{i=1}%
^{n}\underbrace{a_{i,\operatorname*{id}\left(  i\right)  }}_{=a_{i,i}}%
+\sum_{\substack{\sigma\in S_{n};\\\sigma\neq\operatorname*{id}}}\left(
-1\right)  ^{\sigma}\underbrace{\prod_{i=1}^{n}a_{i,\sigma\left(  i\right)  }%
}_{\substack{=0\\\text{(by (\ref{sol.ps4.3.0}))}}}\\
&  \ \ \ \ \ \ \ \ \ \ \left(  \text{here, we have moved the addend for
}\sigma=\operatorname*{id}\text{ out of the sum}\right) \\
&  =\prod_{i=1}^{n}a_{i,i}+\underbrace{\sum_{\substack{\sigma\in
S_{n};\\\sigma\neq\operatorname*{id}}}\left(  -1\right)  ^{\sigma}0}%
_{=0}=\prod_{i=1}^{n}a_{i,i}=a_{1,1}a_{2,2}\cdots a_{n,n}.
\end{align*}
This solves Exercise \ref{exe.ps4.3}.
\end{proof}

\subsection{Solution to Exercise \ref{exe.ps4.4}}

\begin{vershort}
\begin{proof}
[Solution to Exercise \ref{exe.ps4.4}.]The map $S_{n}\rightarrow
S_{n},\ \sigma\mapsto\sigma^{-1}$ (that is, the map from $S_{n}$ to $S_{n}$
which sends every permutation to its inverse) is a bijection\footnote{In fact,
this map is its own inverse: Indeed, every $\sigma\in S_{n}$ satisfies
$\left(  \sigma^{-1}\right)  ^{-1}=\sigma$.}.

Write $A$ in the form $A=\left(  a_{i,j}\right)  _{1\leq i\leq n,\ 1\leq j\leq
n}$. Thus, $A^{T}=\left(  a_{j,i}\right)  _{1\leq i\leq n,\ 1\leq j\leq n}$
(by the definition of $A^{T}$). Hence, (\ref{eq.det.eq.2}) (applied to $A^{T}$
and $a_{j,i}$ instead of $A$ and $a_{i,j}$) yields%
\begin{equation}
\det\left(  A^{T}\right)  =\sum_{\sigma\in S_{n}}\left(  -1\right)  ^{\sigma
}\prod_{i=1}^{n}a_{\sigma\left(  i\right)  ,i}=\sum_{\sigma\in S_{n}}\left(
-1\right)  ^{\sigma}\prod_{i=1}^{n}a_{\sigma^{-1}\left(  i\right)  ,i}
\label{sol.ps4.4.short.3a}%
\end{equation}
(here, we have substituted $\sigma^{-1}$ for $\sigma$ in the sum, since the
map $S_{n}\rightarrow S_{n},\ \sigma\mapsto\sigma^{-1}$ is a bijection). But
every $\sigma\in S_{n}$ satisfies%
\[
\prod_{i=1}^{n}a_{\sigma^{-1}\left(  i\right)  ,i}=\prod_{i=1}^{n}%
a_{i,\sigma\left(  i\right)  }%
\]
\footnote{\textit{Proof.} Let $\sigma\in S_{n}$. Then, $\sigma$ is a
permutation of $\left\{  1,2,\ldots,n\right\}  $, thus a bijection from
$\left\{  1,2,\ldots,n\right\}  $ to $\left\{  1,2,\ldots,n\right\}  $. Hence,
we can substitute $\sigma\left(  i\right)  $ for $i$ in the product
$\prod_{i=1}^{n}a_{\sigma^{-1}\left(  i\right)  ,i}$. We thus obtain%
\[
\prod_{i=1}^{n}a_{\sigma^{-1}\left(  i\right)  ,i}=\prod_{i=1}^{n}%
\underbrace{a_{\sigma^{-1}\left(  \sigma\left(  i\right)  \right)
,\sigma\left(  i\right)  }}_{=a_{i,\sigma\left(  i\right)  }}=\prod_{i=1}%
^{n}a_{i,\sigma\left(  i\right)  },
\]
qed.}. Hence, (\ref{sol.ps4.4.short.3a}) becomes%
\[
\det\left(  A^{T}\right)  =\sum_{\sigma\in S_{n}}\left(  -1\right)  ^{\sigma
}\underbrace{\prod_{i=1}^{n}a_{\sigma^{-1}\left(  i\right)  ,i}}_{=\prod
_{i=1}^{n}a_{i,\sigma\left(  i\right)  }}=\sum_{\sigma\in S_{n}}\left(
-1\right)  ^{\sigma}\prod_{i=1}^{n}a_{i,\sigma\left(  i\right)  }=\det
A\ \ \ \ \ \ \ \ \ \ \left(  \text{by (\ref{eq.det.eq.2})}\right)  .
\]
This solves Exercise \ref{exe.ps4.4}.
\end{proof}
\end{vershort}

\begin{verlong}
\begin{proof}
[Solution to Exercise \ref{exe.ps4.4}.]Define a map $\Phi:S_{n}\rightarrow
S_{n}$ by%
\[
\Phi\left(  \sigma\right)  =\sigma^{-1}\ \ \ \ \ \ \ \ \ \ \text{for every
}\sigma\in S_{n}.
\]
Then, $\Phi\circ\Phi=\operatorname*{id}$\ \ \ \ \footnote{\textit{Proof.}
Every $\sigma\in S_{n}$ satisfies%
\begin{align*}
\left(  \Phi\circ\Phi\right)  \left(  \sigma\right)   &  =\Phi\left(
\underbrace{\Phi\left(  \sigma\right)  }_{=\sigma^{-1}}\right)  =\Phi\left(
\sigma^{-1}\right)  =\left(  \sigma^{-1}\right)  ^{-1}%
\ \ \ \ \ \ \ \ \ \ \left(  \text{by the definition of }\Phi\right) \\
&  =\sigma=\operatorname*{id}\left(  \sigma\right)  .
\end{align*}
Thus, $\Phi\circ\Phi=\operatorname*{id}$, qed.}. Hence, the maps $\Phi$ and
$\Phi$ are mutually inverse. Therefore, the map $\Phi$ is a bijection.

Write $A$ in the form $A=\left(  a_{i,j}\right)  _{1\leq i\leq n,\ 1\leq j\leq
n}$. Thus, $A^{T}=\left(  a_{j,i}\right)  _{1\leq i\leq n,\ 1\leq j\leq n}$
(by the definition of $A^{T}$). Hence, (\ref{eq.det.eq.2}) (applied to $A^{T}$
and $a_{j,i}$ instead of $A$ and $a_{i,j}$) yields%
\begin{align}
\det\left(  A^{T}\right)   &  =\sum_{\sigma\in S_{n}}\left(  -1\right)
^{\sigma}\underbrace{\prod_{i=1}^{n}}_{=\prod_{i\in\left\{  1,2,\ldots
,n\right\}  }}a_{\sigma\left(  i\right)  ,i}=\sum_{\sigma\in S_{n}}\left(
-1\right)  ^{\sigma}\prod_{i\in\left\{  1,2,\ldots,n\right\}  }a_{\sigma
\left(  i\right)  ,i}\nonumber\\
&  =\sum_{\sigma\in S_{n}}\left(  -1\right)  ^{\sigma}\prod_{i\in\left\{
1,2,\ldots,n\right\}  }a_{\left(  \Phi\left(  \sigma\right)  \right)  \left(
i\right)  ,i} \label{sol.ps4.4.3}%
\end{align}
(here, we have substituted $\Phi\left(  \sigma\right)  $ for $\sigma$ in the
sum, since the map $\Phi$ is a bijection). But every $\sigma\in S_{n}$
satisfies%
\[
\prod_{i\in\left\{  1,2,\ldots,n\right\}  }a_{\left(  \Phi\left(
\sigma\right)  \right)  \left(  i\right)  ,i}=\prod_{i\in\left\{
1,2,\ldots,n\right\}  }a_{i,\sigma\left(  i\right)  }%
\]
\footnote{\textit{Proof.} Let $\sigma\in S_{n}$. Then, $\sigma$ is a
permutation of $\left\{  1,2,\ldots,n\right\}  $, thus a bijection from
$\left\{  1,2,\ldots,n\right\}  $ to $\left\{  1,2,\ldots,n\right\}  $. Hence,
we can substitute $\sigma\left(  i\right)  $ for $i$ in the product
$\prod_{i\in\left\{  1,2,\ldots,n\right\}  }a_{\left(  \Phi\left(
\sigma\right)  \right)  \left(  i\right)  ,i}$. We thus obtain%
\begin{align*}
\prod_{i\in\left\{  1,2,\ldots,n\right\}  }a_{\left(  \Phi\left(
\sigma\right)  \right)  \left(  i\right)  ,i}  &  =\prod_{i\in\left\{
1,2,\ldots,n\right\}  }a_{\left(  \Phi\left(  \sigma\right)  \right)  \left(
\sigma\left(  i\right)  \right)  ,\sigma\left(  i\right)  }=\prod
_{i\in\left\{  1,2,\ldots,n\right\}  }\underbrace{a_{\sigma^{-1}\left(
\sigma\left(  i\right)  \right)  ,\sigma\left(  i\right)  }}%
_{\substack{=a_{i,\sigma\left(  i\right)  }\\\text{(since }\sigma^{-1}\left(
\sigma\left(  i\right)  \right)  =i\text{)}}}\ \ \ \ \ \ \ \ \ \ \left(
\text{since }\Phi\left(  \sigma\right)  =\sigma^{-1}\right) \\
&  =\prod_{i\in\left\{  1,2,\ldots,n\right\}  }a_{i,\sigma\left(  i\right)  },
\end{align*}
qed.}. Hence, (\ref{sol.ps4.4.3}) becomes%
\begin{align*}
\det\left(  A^{T}\right)   &  =\sum_{\sigma\in S_{n}}\left(  -1\right)
^{\sigma}\underbrace{\prod_{i\in\left\{  1,2,\ldots,n\right\}  }a_{\left(
\Phi\left(  \sigma\right)  \right)  \left(  i\right)  ,i}}_{=\prod
_{i\in\left\{  1,2,\ldots,n\right\}  }a_{i,\sigma\left(  i\right)  }}%
=\sum_{\sigma\in S_{n}}\left(  -1\right)  ^{\sigma}\underbrace{\prod
_{i\in\left\{  1,2,\ldots,n\right\}  }}_{=\prod_{i=1}^{n}}a_{i,\sigma\left(
i\right)  }\\
&  =\sum_{\sigma\in S_{n}}\left(  -1\right)  ^{\sigma}\prod_{i=1}%
^{n}a_{i,\sigma\left(  i\right)  }=\det A\ \ \ \ \ \ \ \ \ \ \left(  \text{by
(\ref{eq.det.eq.2})}\right)  .
\end{align*}
This solves Exercise \ref{exe.ps4.4}.
\end{proof}
\end{verlong}

\subsection{Solution to Exercise \ref{exe.transpose.basics}}

\begin{proof}
[Solution to Exercise \ref{exe.transpose.basics}.]\textbf{(a)} Let $u$, $v$
and $w$ be three nonnegative integers. Let $P$ be a $u\times v$-matrix, and
let $Q$ be a $v\times w$-matrix. We must prove that $\left(  PQ\right)
^{T}=Q^{T}P^{T}$.

Write the $u\times v$-matrix $P$ in the form $P=\left(  p_{i,j}\right)
_{1\leq i\leq u,\ 1\leq j\leq v}$. Write the $v\times w$-matrix $Q$ in the
form $Q=\left(  q_{i,j}\right)  _{1\leq i\leq v,\ 1\leq j\leq w}$. By the
definition of the product $PQ$, we obtain%
\begin{align*}
PQ  &  =\left(  \sum_{k=1}^{v}\underbrace{p_{i,k}q_{k,j}}_{=q_{k,j}p_{i,k}%
}\right)  _{1\leq i\leq u,\ 1\leq j\leq w}\ \ \ \ \ \ \ \ \ \ \left(
\begin{array}
[c]{c}%
\text{since }P=\left(  p_{i,j}\right)  _{1\leq i\leq u,\ 1\leq j\leq v}\\
\text{and }Q=\left(  q_{i,j}\right)  _{1\leq i\leq v,\ 1\leq j\leq w}%
\end{array}
\right) \\
&  =\left(  \sum_{k=1}^{v}q_{k,j}p_{i,k}\right)  _{1\leq i\leq u,\ 1\leq j\leq
w}.
\end{align*}
Hence,
\begin{align}
\left(  PQ\right)  ^{T}  &  =\left(  \left(  \sum_{k=1}^{v}q_{k,j}%
p_{i,k}\right)  _{1\leq i\leq u,\ 1\leq j\leq w}\right)  ^{T}\nonumber\\
&  =\left(  \sum_{k=1}^{v}q_{k,i}p_{j,k}\right)  _{1\leq i\leq w,\ 1\leq j\leq
u} \label{sol.transpose.basics.a.1}%
\end{align}
(by the definition of the transpose of a matrix).

On the other hand, from $P=\left(  p_{i,j}\right)  _{1\leq i\leq u,\ 1\leq
j\leq v}$, we obtain%
\[
P^{T}=\left(  \left(  p_{i,j}\right)  _{1\leq i\leq u,\ 1\leq j\leq v}\right)
^{T}=\left(  p_{j,i}\right)  _{1\leq i\leq v,\ 1\leq j\leq u}%
\]
(by the definition of the transpose of a matrix).

\begin{vershort}
Likewise, $Q^{T}=\left(  q_{j,i}\right)  _{1\leq i\leq w,\ 1\leq j\leq v}$.
\end{vershort}

\begin{verlong}
Also, from $Q=\left(  q_{i,j}\right)  _{1\leq i\leq v,\ 1\leq j\leq w}$, we
obtain%
\[
Q^{T}=\left(  \left(  q_{i,j}\right)  _{1\leq i\leq v,\ 1\leq j\leq w}\right)
^{T}=\left(  q_{j,i}\right)  _{1\leq i\leq w,\ 1\leq j\leq v}%
\]
(by the definition of the transpose of a matrix).
\end{verlong}

By the definition of the product $Q^{T}P^{T}$, we obtain%
\[
Q^{T}P^{T}=\left(  \sum_{k=1}^{v}q_{k,i}p_{j,k}\right)  _{1\leq i\leq
w,\ 1\leq j\leq u}%
\]
(since $Q^{T}=\left(  q_{j,i}\right)  _{1\leq i\leq w,\ 1\leq j\leq v}$ and
$P^{T}=\left(  p_{j,i}\right)  _{1\leq i\leq v,\ 1\leq j\leq u}$). Comparing
this with (\ref{sol.transpose.basics.a.1}), we obtain $\left(  PQ\right)
^{T}=Q^{T}P^{T}$. This solves Exercise \ref{exe.transpose.basics} \textbf{(a)}.

\textbf{(b)} Let $u\in\mathbb{N}$. We must prove that $\left(  I_{u}\right)
^{T}=I_{u}$.

For any two objects $i$ and $j$, we define an element $\delta_{i,j}%
\in\mathbb{K}$ by $\delta_{i,j}=%
\begin{cases}
1, & \text{if }i=j;\\
0, & \text{if }i\neq j
\end{cases}
$.

\begin{vershort}
Then, it is clear that $\delta_{i,j}=\delta_{j,i}$ for any two objects $i$ and
$j$ (since $i=j$ holds if and only if $j=i$). Thus, $\left(  \delta
_{i,j}\right)  _{1\leq i\leq u,\ 1\leq j\leq u}=\left(  \delta_{j,i}\right)
_{1\leq i\leq u,\ 1\leq j\leq u}$.
\end{vershort}

\begin{verlong}
Then, any two objects $i$ and $j$ satisfy
\begin{equation}
\delta_{i,j}=\delta_{j,i} \label{sol.transpose.basics.b.commut}%
\end{equation}
\footnote{\textit{Proof of (\ref{sol.transpose.basics.b.commut}):} Let $i$ and
$j$ be two objects. The definition of $\delta_{i,j}$ yields
\begin{equation}
\delta_{i,j}=%
\begin{cases}
1, & \text{if }i=j;\\
0, & \text{if }i\neq j
\end{cases}
=%
\begin{cases}
1, & \text{if }j=i;\\
0, & \text{if }j\neq i
\end{cases}
\label{sol.transpose.basics.b.commut.pf.1}%
\end{equation}
(because the condition $\left(  i=j\right)  $ is equivalent to $\left(
j=i\right)  $, and because the condition $\left(  i\neq j\right)  $ is
equivalent to $\left(  j\neq i\right)  $). On the other hand, the definition
of $\delta_{j,i}$ yields%
\[
\delta_{j,i}=%
\begin{cases}
1, & \text{if }j=i;\\
0, & \text{if }j\neq i
\end{cases}
.
\]
Comparing this with (\ref{sol.transpose.basics.b.commut.pf.1}), we obtain
$\delta_{i,j}=\delta_{j,i}$. This proves (\ref{sol.transpose.basics.b.commut}%
).}. Thus,%
\[
\left(  \underbrace{\delta_{i,j}}_{\substack{=\delta_{j,i}\\\text{(by
(\ref{sol.transpose.basics.b.commut}))}}}\right)  _{1\leq i\leq u,\ 1\leq
j\leq u}=\left(  \delta_{j,i}\right)  _{1\leq i\leq u,\ 1\leq j\leq u}.
\]

\end{verlong}

But $I_{u}=\left(  \delta_{i,j}\right)  _{1\leq i\leq u,\ 1\leq j\leq u}$ (by
the definition of $I_{u}$). Hence,
\[
\left(  I_{u}\right)  ^{T}=\left(  \left(  \delta_{i,j}\right)  _{1\leq i\leq
u,\ 1\leq j\leq u}\right)  ^{T}=\left(  \delta_{j,i}\right)  _{1\leq i\leq
u,\ 1\leq j\leq u}%
\]
(by the definition of the transpose of a matrix). Comparing this with%
\[
I_{u}=\left(  \delta_{i,j}\right)  _{1\leq i\leq u,\ 1\leq j\leq u}=\left(
\delta_{j,i}\right)  _{1\leq i\leq u,\ 1\leq j\leq u},
\]
we obtain $\left(  I_{u}\right)  ^{T}=I_{u}$. This solves Exercise
\ref{exe.transpose.basics} \textbf{(b)}.

\textbf{(c)} Let $u$ and $v$ be two nonnegative integers. Let $P$ be a
$u\times v$-matrix. Let $\lambda\in\mathbb{K}$. We must prove that $\left(
\lambda P\right)  ^{T}=\lambda P^{T}$.

Write the $u\times v$-matrix $P$ in the form $P=\left(  p_{i,j}\right)
_{1\leq i\leq u,\ 1\leq j\leq v}$. Thus,
\[
P^{T}=\left(  \left(  p_{i,j}\right)  _{1\leq i\leq u,\ 1\leq j\leq v}\right)
^{T}=\left(  p_{j,i}\right)  _{1\leq i\leq v,\ 1\leq j\leq u}%
\]
(by the definition of the transpose of a matrix). Hence,%
\begin{equation}
\lambda\underbrace{P^{T}}_{=\left(  p_{j,i}\right)  _{1\leq i\leq v,\ 1\leq
j\leq u}}=\lambda\left(  p_{j,i}\right)  _{1\leq i\leq v,\ 1\leq j\leq
u}=\left(  \lambda p_{j,i}\right)  _{1\leq i\leq v,\ 1\leq j\leq u}.
\label{sol.transpose.basics.c.1}%
\end{equation}

On the other hand,
\[
\lambda\underbrace{P}_{=\left(  p_{i,j}\right)  _{1\leq i\leq u,\ 1\leq j\leq
v}}=\lambda\left(  p_{i,j}\right)  _{1\leq i\leq u,\ 1\leq j\leq v}=\left(
\lambda p_{i,j}\right)  _{1\leq i\leq u,\ 1\leq j\leq v}.
\]
Therefore,%
\begin{align*}
\left(  \lambda P\right)  ^{T}  &  =\left(  \left(  \lambda p_{i,j}\right)
_{1\leq i\leq u,\ 1\leq j\leq v}\right)  ^{T}=\left(  \lambda p_{j,i}\right)
_{1\leq i\leq v,\ 1\leq j\leq u}\\
&  \ \ \ \ \ \ \ \ \ \ \ \ \ \ \ \ \ \ \ \ \left(  \text{by the definition of
the transpose of a matrix}\right) \\
&  =\lambda P^{T}%
\end{align*}
(by (\ref{sol.transpose.basics.c.1})). This solves Exercise
\ref{exe.transpose.basics} \textbf{(c)}.

\begin{vershort}
\textbf{(d)} The solution to Exercise \ref{exe.transpose.basics} \textbf{(d)}
is very similar to that of Exercise \ref{exe.transpose.basics} \textbf{(c)},
and we thus omit it.
\end{vershort}

\begin{verlong}
\textbf{(d)} Let $u$ and $v$ be two nonnegative integers. Let $P$ and $Q$ be
two $u\times v$-matrices. We must prove that $\left(  P+Q\right)  ^{T}%
=P^{T}+Q^{T}$.

Write the $u\times v$-matrix $P$ in the form $P=\left(  p_{i,j}\right)
_{1\leq i\leq u,\ 1\leq j\leq v}$. Thus,
\begin{equation}
P^{T}=\left(  \left(  p_{i,j}\right)  _{1\leq i\leq u,\ 1\leq j\leq v}\right)
^{T}=\left(  p_{j,i}\right)  _{1\leq i\leq v,\ 1\leq j\leq u}
\label{sol.transpose.basics.d.P=}%
\end{equation}
(by the definition of the transpose of a matrix).

Write the $u\times v$-matrix $Q$ in the form $Q=\left(  q_{i,j}\right)
_{1\leq i\leq u,\ 1\leq j\leq v}$. Thus,
\begin{equation}
Q^{T}=\left(  \left(  q_{i,j}\right)  _{1\leq i\leq u,\ 1\leq j\leq v}\right)
^{T}=\left(  q_{j,i}\right)  _{1\leq i\leq v,\ 1\leq j\leq u}
\label{sol.transpose.basics.d.Q=}%
\end{equation}
(by the definition of the transpose of a matrix).

Adding together the equalities (\ref{sol.transpose.basics.d.P=}) and
(\ref{sol.transpose.basics.d.Q=}), we obtain%
\begin{align}
P^{T}+Q^{T}  &  =\left(  p_{j,i}\right)  _{1\leq i\leq v,\ 1\leq j\leq
u}+\left(  q_{j,i}\right)  _{1\leq i\leq v,\ 1\leq j\leq u}\nonumber\\
&  =\left(  p_{j,i}+q_{j,i}\right)  _{1\leq i\leq v,\ 1\leq j\leq u}.
\label{sol.transpose.basics.d.2}%
\end{align}

On the other hand, adding together the equalities $P=\left(  p_{i,j}\right)
_{1\leq i\leq u,\ 1\leq j\leq v}$ and $Q=\left(  q_{i,j}\right)  _{1\leq i\leq
u,\ 1\leq j\leq v}$, we obtain%
\[
P+Q=\left(  p_{i,j}\right)  _{1\leq i\leq u,\ 1\leq j\leq v}+\left(
q_{i,j}\right)  _{1\leq i\leq u,\ 1\leq j\leq v}=\left(  p_{i,j}%
+q_{i,j}\right)  _{1\leq i\leq u,\ 1\leq j\leq v}.
\]
Hence,%
\begin{align*}
\left(  P+Q\right)  ^{T}  &  =\left(  \left(  p_{i,j}+q_{i,j}\right)  _{1\leq
i\leq u,\ 1\leq j\leq v}\right)  ^{T}=\left(  p_{j,i}+q_{j,i}\right)  _{1\leq
i\leq v,\ 1\leq j\leq u}\\
&  \ \ \ \ \ \ \ \ \ \ \left(  \text{by the definition of the transpose of a
matrix}\right) \\
&  =P^{T}+Q^{T}\ \ \ \ \ \ \ \ \ \ \left(  \text{by
(\ref{sol.transpose.basics.d.2})}\right)  .
\end{align*}
This solves Exercise \ref{exe.transpose.basics} \textbf{(d)}.
\end{verlong}

\textbf{(e)} Let $u$ and $v$ be two nonnegative integers. Let $P$ be a
$u\times v$-matrix. We must prove that $\left(  P^{T}\right)  ^{T}=P$.

Write the $u\times v$-matrix $P$ in the form $P=\left(  p_{i,j}\right)
_{1\leq i\leq u,\ 1\leq j\leq v}$. Thus,
\[
P^{T}=\left(  \left(  p_{i,j}\right)  _{1\leq i\leq u,\ 1\leq j\leq v}\right)
^{T}=\left(  p_{j,i}\right)  _{1\leq i\leq v,\ 1\leq j\leq u}%
\]
(by the definition of the transpose of a matrix). Hence,%
\begin{align*}
\left(  P^{T}\right)  ^{T}  &  =\left(  \left(  p_{j,i}\right)  _{1\leq i\leq
v,\ 1\leq j\leq u}\right)  ^{T}=\left(  p_{i,j}\right)  _{1\leq i\leq
u,\ 1\leq j\leq v}\\
&  \ \ \ \ \ \ \ \ \ \ \left(  \text{by the definition of the transpose of a
matrix}\right) \\
&  =P.
\end{align*}
This solves Exercise \ref{exe.transpose.basics} \textbf{(e)}.
\end{proof}

\subsection{Solution to Exercise \ref{exe.ps4.5}}

\begin{proof}
[Solution to Exercise \ref{exe.ps4.5}.]Of course, both parts of Exercise
\ref{exe.ps4.5} can be solved directly using (\ref{eq.det.eq.1}). This
solution, however, is tedious (particularly for part \textbf{(b)} of this
exercise). Let us show a smarter way.

\textbf{(a)} Let $A$ be the matrix $\left(
\begin{array}
[c]{cccc}%
a & b & c & d\\
\ell & 0 & 0 & e\\
k & 0 & 0 & f\\
j & i & h & g
\end{array}
\right)  $. We want to find $\det A$.

We write the matrix $A$ in the form $A=\left(  a_{u,v}\right)  _{1\leq
u\leq4,\ 1\leq v\leq4}$\ \ \ \ \footnote{We cannot write $A=\left(
a_{i,j}\right)  _{1\leq i\leq4,\ 1\leq j\leq4}$ because the letters $i$ and
$j$ are already taken for something different.}. Thus,%
\begin{align*}
a_{1,1}  &  =a,\ \ \ \ \ \ \ \ \ \ a_{1,2}=b,\ \ \ \ \ \ \ \ \ \ a_{1,3}%
=c,\ \ \ \ \ \ \ \ \ \ a_{1,4}=d,\\
a_{2,1}  &  =\ell,\ \ \ \ \ \ \ \ \ \ a_{2,2}=0,\ \ \ \ \ \ \ \ \ \ a_{2,3}%
=0,\ \ \ \ \ \ \ \ \ \ a_{2,4}=e,\\
a_{3,1}  &  =k,\ \ \ \ \ \ \ \ \ \ a_{3,2}=0,\ \ \ \ \ \ \ \ \ \ a_{3,3}%
=0,\ \ \ \ \ \ \ \ \ \ a_{3,4}=f,\\
a_{4,1}  &  =j,\ \ \ \ \ \ \ \ \ \ a_{4,2}=i,\ \ \ \ \ \ \ \ \ \ a_{4,3}%
=h,\ \ \ \ \ \ \ \ \ \ a_{4,4}=g.
\end{align*}
Now, applying (\ref{eq.det.eq.1}) to $n=4$, we obtain%
\begin{equation}
\det A=\sum_{\sigma\in S_{4}}\left(  -1\right)  ^{\sigma}a_{1,\sigma\left(
1\right)  }a_{2,\sigma\left(  2\right)  }a_{3,\sigma\left(  3\right)
}a_{4,\sigma\left(  4\right)  }. \label{sol.ps4.5.a.1}%
\end{equation}

The sum on the right hand side of (\ref{sol.ps4.5.a.1}) has $\left\vert
S_{4}\right\vert =4!=24$ addends. However, some of them are $0$. Namely, every
addend corresponding to a permutation $\sigma\in S_{4}$ satisfying
$\sigma\left(  2\right)  \notin\left\{  1,4\right\}  $ must be $0$%
\ \ \ \ \footnote{\textit{Proof.} Let $\sigma\in S_{4}$ be such that
$\sigma\left(  2\right)  \notin\left\{  1,4\right\}  $. We must then show that
the addend on the right hand side of (\ref{sol.ps4.5.a.1}) corresponding to
this $\sigma$ must be $0$. In other words, we have to show that $\left(
-1\right)  ^{\sigma}a_{1,\sigma\left(  1\right)  }a_{2,\sigma\left(  2\right)
}a_{3,\sigma\left(  3\right)  }a_{4,\sigma\left(  4\right)  }=0$.
\par
We have $\sigma\left(  2\right)  \notin\left\{  1,4\right\}  $, and thus
$\sigma\left(  2\right)  \in\left\{  2,3\right\}  $. Hence, $a_{2,\sigma
\left(  2\right)  }=0$ (because $a_{2,2}=0$ and $a_{2,3}=0$), and thus
$\left(  -1\right)  ^{\sigma}a_{1,\sigma\left(  1\right)  }%
\underbrace{a_{2,\sigma\left(  2\right)  }}_{=0}a_{3,\sigma\left(  3\right)
}a_{4,\sigma\left(  4\right)  }=0$, qed.}. Hence, all such addends can be
removed from the sum (without changing the value of this sum). Similarly, all
addends corresponding to permutations $\sigma\in S_{4}$ satisfying
$\sigma\left(  3\right)  \notin\left\{  1,4\right\}  $ must be $0$, and can
therefore also be removed from the sum. The addends that survive these two
removals are the ones that correspond to permutations $\sigma\in S_{4}$
satisfying $\sigma\left(  2\right)  \in\left\{  1,4\right\}  $ and
$\sigma\left(  3\right)  \in\left\{  1,4\right\}  $. It is easy to see that
there are exactly four such permutations: In one-line notation, these
permutations are $\left(  2,1,4,3\right)  $, $\left(  2,4,1,3\right)  $,
$\left(  3,1,4,2\right)  $ and $\left(  3,4,1,2\right)  $. The addends
corresponding to these permutations are $a_{1,2}a_{2,1}a_{3,4}a_{4,3}$,
$-a_{1,2}a_{2,4}a_{3,1}a_{4,3}$, $-a_{1,3}a_{2,1}a_{3,4}a_{4,2}$ and
$a_{1,3}a_{2,4}a_{3,1}a_{4,2}$. Hence, (\ref{sol.ps4.5.a.1}) simplifies to%
\begin{align*}
&  \det A\\
&  =\underbrace{a_{1,2}}_{=b}\underbrace{a_{2,1}}_{=\ell}\underbrace{a_{3,4}%
}_{=f}\underbrace{a_{4,3}}_{=h}-\underbrace{a_{1,2}}_{=b}\underbrace{a_{2,4}%
}_{=e}\underbrace{a_{3,1}}_{=k}\underbrace{a_{4,3}}_{=h}-\underbrace{a_{1,3}%
}_{=c}\underbrace{a_{2,1}}_{=\ell}\underbrace{a_{3,4}}_{=f}\underbrace{a_{4,2}%
}_{=i}+\underbrace{a_{1,3}}_{=c}\underbrace{a_{2,4}}_{=e}\underbrace{a_{3,1}%
}_{=k}\underbrace{a_{4,2}}_{=i}\\
&  =b\ell fh-bekh-c\ell fi+ceki.
\end{align*}
This is a simple enough formula to consider an answer to Exercise
\ref{exe.ps4.5} \textbf{(a)}, but we can simplify it even further. Namely,%
\[
\det A=b\ell fh-bekh-c\ell fi+ceki=\left(  bh-ci\right)  \left(  \ell
f-ek\right)  .
\]
Exercise \ref{exe.ps4.5} \textbf{(a)} is solved.\footnote{It is not an
accident that this determinant factors so nicely! See Example
\ref{exa.det.blocktria-twisted.exas} \textbf{(c)} for a more conceptual way to
see this.}

\textbf{(b)} Let $A$ be the matrix $\left(
\begin{array}
[c]{ccccc}%
a & b & c & d & e\\
f & 0 & 0 & 0 & g\\
h & 0 & 0 & 0 & i\\
j & 0 & 0 & 0 & k\\
\ell & m & n & o & p
\end{array}
\right)  $. We want to find $\det A$.

We write the matrix $A$ in the form
\[
A=\left(  a_{u,v}\right)  _{1\leq u\leq5,\ 1\leq v\leq5}.
\]
Thus, $a_{1,1}=a$, $a_{1,2}=b$, etc.. For us, the most important property of
$A$ is that the $3\times3$-submatrix in the middle of $A$ is filled with
zeroes. In other words,%
\begin{equation}
a_{u,v}=0\ \ \ \ \ \ \ \ \ \ \text{for every }u\in\left\{  2,3,4\right\}
\text{ and }v\in\left\{  2,3,4\right\}  . \label{sol.ps4.5.b.2}%
\end{equation}

Now, applying (\ref{eq.det.eq.1}) to $n=5$, we obtain%
\begin{equation}
\det A=\sum_{\sigma\in S_{5}}\left(  -1\right)  ^{\sigma}a_{1,\sigma\left(
1\right)  }a_{2,\sigma\left(  2\right)  }a_{3,\sigma\left(  3\right)
}a_{4,\sigma\left(  4\right)  }a_{5,\sigma\left(  5\right)  }.
\label{sol.ps4.5.b.3}%
\end{equation}
But every $\sigma\in S_{5}$ satisfies $a_{2,\sigma\left(  2\right)
}a_{3,\sigma\left(  3\right)  }a_{4,\sigma\left(  4\right)  }=0$%
\ \ \ \ \footnote{\textit{Proof.} Let $\sigma\in S_{5}$. Then, $\sigma$ is a
permutation of $\left\{  1,2,3,4,5\right\}  $, and thus an injective map.
Therefore, the numbers $\sigma\left(  2\right)  ,\sigma\left(  3\right)
,\sigma\left(  4\right)  $ are pairwise distinct.
\par
We now claim that there exists an $u\in\left\{  2,3,4\right\}  $ such that
$\sigma\left(  u\right)  \in\left\{  2,3,4\right\}  $. In order to prove this,
we assume the contrary. Thus, every $u\in\left\{  2,3,4\right\}  $ satisfies
$\sigma\left(  u\right)  \notin\left\{  2,3,4\right\}  $. Hence, every
$u\in\left\{  2,3,4\right\}  $ satisfies $\sigma\left(  u\right)  \in\left\{
1,5\right\}  $ (since $\sigma\left(  u\right)  \in\left\{  1,2,3,4,5\right\}
$ but $\sigma\left(  u\right)  \notin\left\{  2,3,4\right\}  $). In other
words, the numbers $\sigma\left(  2\right)  ,\sigma\left(  3\right)
,\sigma\left(  4\right)  $ belong to $\left\{  1,5\right\}  $. Hence,
$\sigma\left(  2\right)  ,\sigma\left(  3\right)  ,\sigma\left(  4\right)  $
are three distinct numbers belonging to the set $\left\{  1,5\right\}  $. But
this is absurd, since the set $\left\{  1,5\right\}  $ does not contain three
distinct numbers. Hence, we have obtained a contradiction. This shows that our
assumption was wrong.
\par
We thus have shown that there exists an $u\in\left\{  2,3,4\right\}  $ such
that $\sigma\left(  u\right)  \in\left\{  2,3,4\right\}  $. Consider such a
$u$. Applying (\ref{sol.ps4.5.b.2}) to $v=\sigma\left(  u\right)  $, we now
obtain $a_{u,\sigma\left(  u\right)  }=0$. But $u\in\left\{  2,3,4\right\}  $,
so that $a_{u,\sigma\left(  u\right)  }$ is a factor in the product
$a_{2,\sigma\left(  2\right)  }a_{3,\sigma\left(  3\right)  }a_{4,\sigma
\left(  4\right)  }$. Hence, the product $a_{2,\sigma\left(  2\right)
}a_{3,\sigma\left(  3\right)  }a_{4,\sigma\left(  4\right)  }$ is $0$ (since
its factor $a_{u,\sigma\left(  u\right)  }$ is $0$), qed.}. Hence,
(\ref{sol.ps4.5.b.3}) becomes%
\[
\det A=\sum_{\sigma\in S_{5}}\left(  -1\right)  ^{\sigma}a_{1,\sigma\left(
1\right)  }\underbrace{a_{2,\sigma\left(  2\right)  }a_{3,\sigma\left(
3\right)  }a_{4,\sigma\left(  4\right)  }}_{=0}a_{5,\sigma\left(  5\right)
}=\sum_{\sigma\in S_{5}}\left(  -1\right)  ^{\sigma}a_{1,\sigma\left(
1\right)  }0a_{5,\sigma\left(  5\right)  }=0.
\]
Exercise \ref{exe.ps4.5} \textbf{(b)} is thus solved.
\end{proof}

\subsection{Solution to Exercise \ref{exe.ps4.6}}

Our solution to Exercise \ref{exe.ps4.6} relies on Lemma \ref{lem.det.sigma}.
Thus, the reader is advised to read the proof of said lemma before the
solution. Furthermore, we shall use the following simple fact:

\begin{lemma}
\label{lem.colsA=rowsATT}Let $n\in\mathbb{N}$ and $m\in\mathbb{N}$. Let $A$ be
an $n\times m$-matrix. Then, the columns of $A$ are the transposes of the
respective rows of $A^{T}$.
\end{lemma}

\begin{proof}
[Proof of Lemma \ref{lem.colsA=rowsATT}.]Fix $k\in\left\{  1,2,\ldots
,m\right\}  $.

Write the $n\times m$-matrix $A$ in the form $A=\left(  a_{i,j}\right)
_{1\leq i\leq n,\ 1\leq j\leq m}$. Then, $A^{T}=\left(  a_{j,i}\right)
_{1\leq i\leq m,\ 1\leq j\leq n}$ (by the definition of $A^{T}$). Hence,%
\[
\left(  \text{the }k\text{-th row of }A^{T}\right)  =\left(  a_{j,k}\right)
_{1\leq i\leq1,\ 1\leq j\leq n}.
\]
Thus,%
\begin{align}
&  \left(  \text{the transpose of }\underbrace{\text{the }k\text{-th row of
}A^{T}}_{=\left(  a_{j,k}\right)  _{1\leq i\leq1,\ 1\leq j\leq n}}\right)
\nonumber\\
&  =\left(  \text{the transpose of }\left(  a_{j,k}\right)  _{1\leq
i\leq1,\ 1\leq j\leq n}\right)  =\left(  \left(  a_{j,k}\right)  _{1\leq
i\leq1,\ 1\leq j\leq n}\right)  ^{T}\nonumber\\
&  =\left(  a_{i,k}\right)  _{1\leq i\leq n,\ 1\leq j\leq1}
\label{pf.lem.colsA=rowsATT.1}%
\end{align}
(by the definition of $\left(  \left(  a_{j,k}\right)  _{1\leq i\leq1,\ 1\leq
j\leq n}\right)  ^{T}$). On the other hand, we have $A=\left(  a_{i,j}\right)
_{1\leq i\leq n,\ 1\leq j\leq m}$, and thus%
\begin{align}
\left(  \text{the }k\text{-th column of }A\right)   &  =\left(  a_{i,k}%
\right)  _{1\leq i\leq n,\ 1\leq j\leq1}\nonumber\\
&  =\left(  \text{the transpose of the }k\text{-th row of }A^{T}\right)
\label{pf.lem.colsA=rowsATT.claim}%
\end{align}
(by (\ref{pf.lem.colsA=rowsATT.1})).

Now, forget that we fixed $k$. We thus have proven that
(\ref{pf.lem.colsA=rowsATT.claim}) holds for each $k\in\left\{  1,2,\ldots
,m\right\}  $. In other words, for each $k\in\left\{  1,2,\ldots,m\right\}  $,
the $k$-th column of $A$ is the transpose of the $k$-th row of $A^{T}$. In
other words, the columns of $A$ are the transposes of the respective rows of
$A^{T}$. This proves Lemma \ref{lem.colsA=rowsATT}.
\end{proof}

\begin{remark}
\label{rmk.colsA=rowsATT.rmk}Let $n\in\mathbb{N}$. Let $A$ be an $n\times
n$-matrix. Lemma \ref{lem.colsA=rowsATT} (applied to $m=n$) shows that the
columns of $A$ are the transposes of the respective rows of $A^{T}$. Thus, we
have a correspondence between the columns of $A$ and the rows of $A^{T}$. We
can use this correspondence to \textquotedblleft transport\textquotedblright%
\ information about the columns of $A$ to the rows of $A^{T}$ and vice versa;
for example:

\begin{itemize}
\item If a column of $A$ consists of zeroes, then the corresponding row of
$A^{T}$ consists of zeroes. The converse also holds.

\item If two columns of $A$ are equal, then the corresponding two rows of
$A^{T}$ are equal. The converse also holds.

\item Let $k\in\left\{  1,2,\ldots,n\right\}  $ and $\lambda\in\mathbb{K}$. If
$B$ is the $n\times n$-matrix obtained from $A$ by multiplying the $k$-th
column by $\lambda$, then $B^{T}$ is the $n\times n$-matrix obtained from
$A^{T}$ by multiplying the $k$-th row by $\lambda$. Again, the converse also holds.

\item Let $k\in\left\{  1,2,\ldots,n\right\}  $. Let $A^{\prime}$ be an
$n\times n$-matrix whose columns equal the corresponding columns of $A$ except
(perhaps) the $k$-th column. Then, $\left(  A^{\prime}\right)  ^{T}$ is an
$n\times n$-matrix whose rows equal the corresponding rows of $A^{T}$ except
(perhaps) the $k$-th row. Again, the converse also holds.

\item Let $k\in\left\{  1,2,\ldots,n\right\}  $. Let $A^{\prime}$ be a further
$n\times n$-matrix. Let $B$ be the $n\times n$-matrix obtained from $A$ by
adding the $k$-th column of $A^{\prime}$ to the $k$-th column of $A$. Then,
$B^{T}$ is the $n\times n$-matrix obtained from $A^{T}$ by adding the $k$-th
row of $\left(  A^{\prime}\right)  ^{T}$ to the $k$-th row of $A^{T}$. Again,
the converse holds as well.
\end{itemize}
\end{remark}

\begin{proof}
[Solution to Exercise \ref{exe.ps4.6}.]Let us write the matrix $A$ in the form
$A=\left(  a_{i,j}\right)  _{1\leq i\leq n,\ 1\leq j\leq n}$. Thus, every
$\left(  i,j\right)  \in\left\{  1,2,\ldots,n\right\}  ^{2}$ satisfies%
\begin{equation}
\left(  \text{the }\left(  i,j\right)  \text{-th entry of the matrix
}A\right)  =a_{i,j}. \label{sol.ps4.6.aij}%
\end{equation}

We let $\left[  n\right]  $ denote the set $\left\{  1,2,\ldots,n\right\}  $.

\textbf{(a)} Let $B$ be an $n\times n$-matrix obtained from $A$ by swapping
two rows. Thus, there exist two distinct elements $u$ and $v$ of $\left\{
1,2,\ldots,n\right\}  $ such that $B$ is the $n\times n$-matrix obtained from
$A$ by swapping the $u$-th row with the $v$-th row. Consider these $u$ and $v$.

We write the matrix $B$ in the form $B=\left(  b_{i,j}\right)  _{1\leq i\leq
n,\ 1\leq j\leq n}$.

\begin{vershort}
Consider the transposition $t_{u,v}$ in $S_{n}$ (defined according to
Definition \ref{def.transpos}). Clearly, $t_{u,v}$ is a permutation of the set
$\left\{  1,2,\ldots,n\right\}  =\left[  n\right]  $, thus a map $\left[
n\right]  \rightarrow\left[  n\right]  $. Also, $\left(  -1\right)  ^{t_{u,v}%
}=-1$ (by Exercise \ref{exe.ps4.1ab} \textbf{(b)}, applied to $i=u$ and $j=v$).

Recall that $B$ is the $n\times n$-matrix obtained from $A$ by swapping the
$u$-th row with the $v$-th row. In other words, for all $k\in\left\{
1,2,\ldots,n\right\}  $, we have
\begin{equation}
\left(  \text{the }k\text{-th row of }B\right)  =\left(  \text{the }%
t_{u,v}\left(  k\right)  \text{-th row of }A\right)
\label{sol.ps4.6.a.short.krow2}%
\end{equation}
(because $t_{u,v}$ is the permutation of $\left\{  1,2,\ldots,n\right\}  $
that swaps $u$ with $v$ while leaving all other numbers fixed). Therefore,
every $\left(  i,j\right)  \in\left\{  1,2,\ldots,n\right\}  ^{2}$ satisfies%
\[
b_{i,j}=a_{t_{u,v}\left(  i\right)  ,j}%
\]
\footnote{\textit{Proof.} We have $B=\left(  b_{i,j}\right)  _{1\leq i\leq
n,\ 1\leq j\leq n}$. Thus, every $\left(  i,j\right)  \in\left\{
1,2,\ldots,n\right\}  ^{2}$ satisfies%
\begin{equation}
\left(  \text{the }\left(  i,j\right)  \text{-th entry of the matrix
}B\right)  =b_{i,j}. \label{sol.ps4.6.a.short.bij.pf.1}%
\end{equation}
Now, let $\left(  i,j\right)  \in\left\{  1,2,\ldots,n\right\}  ^{2}$. Then,
\begin{align*}
b_{i,j}  &  =\left(  \text{the }\left(  i,j\right)  \text{-th entry of the
matrix }B\right)  \ \ \ \ \ \ \ \ \ \ \left(  \text{by
(\ref{sol.ps4.6.a.short.bij.pf.1})}\right) \\
&  =\left(  \text{the }j\text{-th entry of }\underbrace{\text{the }i\text{-th
row of }B}_{\substack{=\left(  \text{the }t_{u,v}\left(  i\right)  \text{-th
row of }A\right)  \\\text{(by (\ref{sol.ps4.6.a.short.krow2}), applied to
}k=i\text{)}}}\right) \\
&  =\left(  \text{the }j\text{-th entry of the }t_{u,v}\left(  i\right)
\text{-th row of }A\right) \\
&  =a_{t_{u,v}\left(  i\right)  ,j}\ \ \ \ \ \ \ \ \ \ \left(  \text{by
(\ref{sol.ps4.6.aij}), applied to }t_{u,v}\left(  i\right)  \text{ instead of
}i\right)  ,
\end{align*}
qed.}. Hence, $B=\left(  \underbrace{b_{i,j}}_{=a_{t_{u,v}\left(  i\right)
,j}}\right)  _{1\leq i\leq n,\ 1\leq j\leq n}=\left(  a_{t_{u,v}\left(
i\right)  ,j}\right)  _{1\leq i\leq n,\ 1\leq j\leq n}$. Therefore, we can
apply Lemma \ref{lem.det.sigma} \textbf{(a)} to $t_{u,v}$, $A$, $a_{i,j}$ and
$B$ instead of $\kappa$, $B$, $b_{i,j}$ and $B_{\kappa}$. We thus obtain $\det
B=\underbrace{\left(  -1\right)  ^{t_{u,v}}}_{=-1}\cdot\det A=-\det A$.
Exercise \ref{exe.ps4.6} \textbf{(a)} is thus solved.
\end{vershort}

\begin{verlong}
We know that $B$ is the $n\times n$-matrix obtained from $A$ by swapping the
$u$-th row with the $v$-th row. In other words,%
\begin{align*}
\left(  \text{the }u\text{-th row of }B\right)   &  =\left(  \text{the
}v\text{-th row of }A\right)  ,\\
\left(  \text{the }v\text{-th row of }B\right)   &  =\left(  \text{the
}u\text{-th row of }A\right)  ,
\end{align*}
and%
\begin{align}
&  \left(  \left(  \text{the }k\text{-th row of }B\right)  =\left(  \text{the
}k\text{-th row of }A\right)  \right. \label{sol.ps4.6.a.krow}\\
&  \ \ \ \ \ \ \ \ \ \ \left.  \text{for all }k\in\left\{  1,2,\ldots
,n\right\}  \text{ satisfying }k\notin\left\{  u,v\right\}  \right)
.\nonumber
\end{align}

Consider the transposition $t_{u,v}$ in $S_{n}$ (defined according to
Definition \ref{def.transpos}). Clearly, $t_{u,v}$ is a permutation of the set
$\left\{  1,2,\ldots,n\right\}  =\left[  n\right]  $, thus a map $\left[
n\right]  \rightarrow\left[  n\right]  $. Also, $\left(  -1\right)  ^{t_{u,v}%
}=-1$ (by Exercise \ref{exe.ps4.1ab} \textbf{(b)}, applied to $i=u$ and $j=v$).

The permutation $t_{u,v}$ is the permutation in $S_{n}$ which swaps $u$ with
$v$ while leaving all other elements of $\left\{  1,2,\ldots,n\right\}  $
unchanged (according to the definition of $t_{u,v}$). In other words, we have
$t_{u,v}\left(  u\right)  =v$, $t_{u,v}\left(  v\right)  =u$ and%
\begin{equation}
\left(  t_{u,v}\left(  k\right)  =k\ \ \ \ \ \ \ \ \ \ \text{for all }%
k\in\left\{  1,2,\ldots,n\right\}  \text{ satisfying }k\notin\left\{
u,v\right\}  \right)  . \label{sol.ps4.6.a.tuv}%
\end{equation}

Now, for all $k\in\left\{  1,2,\ldots,n\right\}  $, we have
\begin{equation}
\left(  \text{the }k\text{-th row of }B\right)  =\left(  \text{the }%
t_{u,v}\left(  k\right)  \text{-th row of }A\right)  \label{sol.ps4.6.a.krow2}%
\end{equation}
\footnote{\textit{Proof of (\ref{sol.ps4.6.a.krow2}):} Let $k\in\left\{
1,2,\ldots,n\right\}  $. We need to prove (\ref{sol.ps4.6.a.krow2}). We are in
one of the following three cases:
\par
\textit{Case 1:} We have $k=u$.
\par
\textit{Case 2:} We have $k=v$.
\par
\textit{Case 3:} We have $k\notin\left\{  u,v\right\}  $.
\par
Let us first consider Case 1. In this case, we have $k=u$. Thus,
$t_{u,v}\left(  \underbrace{k}_{=u}\right)  =t_{u,v}\left(  u\right)  =v$, so
that $\left(  \text{the }t_{u,v}\left(  k\right)  \text{-th row of }A\right)
=\left(  \text{the }v\text{-th row of }A\right)  $. On the other hand,%
\begin{align*}
\left(  \text{the }\underbrace{k}_{=u}\text{-th row of }B\right)   &  =\left(
\text{the }u\text{-th row of }B\right)  =\left(  \text{the }v\text{-th row of
}A\right) \\
&  =\left(  \text{the }t_{u,v}\left(  k\right)  \text{-th row of }A\right)  .
\end{align*}
Thus, (\ref{sol.ps4.6.a.krow2}) is proven in Case 1.
\par
Let us now consider Case 2. In this case, we have $k=v$. Thus, $t_{u,v}\left(
\underbrace{k}_{=v}\right)  =t_{u,v}\left(  v\right)  =u$, so that $\left(
\text{the }t_{u,v}\left(  k\right)  \text{-th row of }A\right)  =\left(
\text{the }u\text{-th row of }A\right)  $. On the other hand,%
\begin{align*}
\left(  \text{the }\underbrace{k}_{=v}\text{-th row of }B\right)   &  =\left(
\text{the }v\text{-th row of }B\right)  =\left(  \text{the }u\text{-th row of
}A\right) \\
&  =\left(  \text{the }t_{u,v}\left(  k\right)  \text{-th row of }A\right)  .
\end{align*}
Thus, (\ref{sol.ps4.6.a.krow2}) is proven in Case 2.
\par
Finally, let us consider Case 3. In this case, we have $k\notin\left\{
u,v\right\}  $. Thus, $t_{u,v}\left(  k\right)  =k$ (by (\ref{sol.ps4.6.a.tuv}%
)), so that $\left(  \text{the }t_{u,v}\left(  k\right)  \text{-th row of
}A\right)  =\left(  \text{the }k\text{-th row of }A\right)  $. On the other
hand,%
\begin{align*}
\left(  \text{the }k\text{-th row of }B\right)   &  =\left(  \text{the
}k\text{-th row of }A\right)  \ \ \ \ \ \ \ \ \ \ \left(  \text{by
(\ref{sol.ps4.6.a.krow})}\right) \\
&  =\left(  \text{the }t_{u,v}\left(  k\right)  \text{-th row of }A\right)  .
\end{align*}
Thus, (\ref{sol.ps4.6.a.krow2}) is proven in Case 3.
\par
We have now proven (\ref{sol.ps4.6.a.krow2}) in each of the three Cases 1, 2
and 3. Hence, (\ref{sol.ps4.6.a.krow2}) always holds, qed.}. Therefore, every
$\left(  i,j\right)  \in\left\{  1,2,\ldots,n\right\}  ^{2}$ satisfies%
\begin{equation}
b_{i,j}=a_{t_{u,v}\left(  i\right)  ,j} \label{sol.ps4.6.a.bij}%
\end{equation}
\footnote{\textit{Proof of (\ref{sol.ps4.6.a.bij}):} We have $B=\left(
b_{i,j}\right)  _{1\leq i\leq n,\ 1\leq j\leq n}$. Thus, every $\left(
i,j\right)  \in\left\{  1,2,\ldots,n\right\}  ^{2}$ satisfies%
\begin{equation}
\left(  \text{the }\left(  i,j\right)  \text{-th entry of the matrix
}B\right)  =b_{i,j}. \label{sol.ps4.6.a.bij.pf.1}%
\end{equation}
Now, let $\left(  i,j\right)  \in\left\{  1,2,\ldots,n\right\}  ^{2}$. Then,
\begin{align*}
b_{i,j}  &  =\left(  \text{the }\left(  i,j\right)  \text{-th entry of the
matrix }B\right)  \ \ \ \ \ \ \ \ \ \ \left(  \text{by
(\ref{sol.ps4.6.a.bij.pf.1})}\right) \\
&  =\left(  \text{the }j\text{-th entry of }\underbrace{\text{the }i\text{-th
row of }B}_{\substack{=\left(  \text{the }t_{u,v}\left(  i\right)  \text{-th
row of }A\right)  \\\text{(by (\ref{sol.ps4.6.a.krow2}), applied to
}k=i\text{)}}}\right) \\
&  =\left(  \text{the }j\text{-th entry of the }t_{u,v}\left(  i\right)
\text{-th row of }A\right) \\
&  =a_{t_{u,v}\left(  i\right)  ,j}\ \ \ \ \ \ \ \ \ \ \left(  \text{by
(\ref{sol.ps4.6.aij}), applied to }t_{u,v}\left(  i\right)  \text{ instead of
}i\right)  ,
\end{align*}
qed.}. Hence, $B=\left(  \underbrace{b_{i,j}}_{=a_{t_{u,v}\left(  i\right)
,j}}\right)  _{1\leq i\leq n,\ 1\leq j\leq n}=\left(  a_{t_{u,v}\left(
i\right)  ,j}\right)  _{1\leq i\leq n,\ 1\leq j\leq n}$. Therefore, we can
apply Lemma \ref{lem.det.sigma} \textbf{(a)} to $t_{u,v}$, $A$, $a_{i,j}$ and
$B$ instead of $\kappa$, $B$, $b_{i,j}$ and $B_{\kappa}$. We thus obtain $\det
B=\underbrace{\left(  -1\right)  ^{t_{u,v}}}_{=-1}\cdot\det A=-\det A$.
Exercise \ref{exe.ps4.6} \textbf{(a)} is thus solved.
\end{verlong}

\textbf{(b)} Let $B$ be an $n\times n$-matrix obtained from $A$ by swapping
two columns. Thus, $B^{T}$ is an $n\times n$-matrix obtained from $A^{T}$ by
swapping two rows (because the columns of $A$ correspond to the rows of
$A^{T}$\ \ \ \ \footnote{See Remark \ref{rmk.colsA=rowsATT.rmk} for the
meaning of \textquotedblleft correspond\textquotedblright\ we are using
here.}). Hence, Exercise \ref{exe.ps4.6} \textbf{(a)} (applied to $A^{T}$ and
$B^{T}$ instead of $A$ and $B$) yields $\det\left(  B^{T}\right)
=-\det\left(  A^{T}\right)  $.

But Exercise \ref{exe.ps4.4} yields $\det\left(  A^{T}\right)  =\det A$. Also,
Exercise \ref{exe.ps4.4} (applied to $B$ instead of $A$) yields $\det\left(
B^{T}\right)  =\det B$. But recall that $\det\left(  B^{T}\right)
=-\det\left(  A^{T}\right)  $. This rewrites as $\det B=-\det A$ (since
$\det\left(  B^{T}\right)  =\det B$ and $\det\left(  A^{T}\right)  =\det A$).
This solves Exercise \ref{exe.ps4.6} \textbf{(b)}.

\textbf{(c)} Assume that a row of $A$ consists of zeroes. Thus, there exists a
$u\in\left\{  1,2,\ldots,n\right\}  $ such that the $u$-th row of $A$ consists
of zeroes. Consider this $u$.

\begin{vershort}
The $u$-th row of the matrix $A$ consists of zeroes. In other words,
\begin{equation}
a_{u,j}=0\ \ \ \ \ \ \ \ \ \ \text{for every }j\in\left\{  1,2,\ldots
,n\right\}  . \label{sol.ps4.6.c.short.1}%
\end{equation}
Now, every $\sigma\in S_{n}$ satisfies $\prod_{i=1}^{n}a_{i,\sigma\left(
i\right)  }=0$\ \ \ \ \footnote{\textit{Proof.} Let $\sigma\in S_{n}$. Then,
the $u$-th factor of the product $\prod_{i=1}^{n}a_{i,\sigma\left(  i\right)
}$ is $a_{u,\sigma\left(  u\right)  }=0$ (by (\ref{sol.ps4.6.c.short.1}),
applied to $j=\sigma\left(  u\right)  $). Hence, the whole product is $0$. In
other words, we have $\prod_{i=1}^{n}a_{i,\sigma\left(  i\right)  }=0$, qed.}.
Thus, (\ref{eq.det.eq.2}) shows that%
\[
\det A=\sum_{\sigma\in S_{n}}\left(  -1\right)  ^{\sigma}\underbrace{\prod
_{i=1}^{n}a_{i,\sigma\left(  i\right)  }}_{=0}=\sum_{\sigma\in S_{n}}\left(
-1\right)  ^{\sigma}0=0.
\]
This solves Exercise \ref{exe.ps4.6} \textbf{(c)}.
\end{vershort}

\begin{verlong}
We have $\left(  \text{the }u\text{-th row of }A\right)  =\underbrace{\left(
0,0,\ldots,0\right)  }_{n\text{ zeroes}}$ (since the $u$-th row of $A$
consists of zeroes). Thus,%
\begin{equation}
a_{u,j}=0\ \ \ \ \ \ \ \ \ \ \text{for every }j\in\left\{  1,2,\ldots
,n\right\}  \label{sol.ps4.6.c.1}%
\end{equation}
\footnote{\textit{Proof of (\ref{sol.ps4.6.c.1}):} Let every $j\in\left\{
1,2,\ldots,n\right\}  $. Applying (\ref{sol.ps4.6.aij}) to $i=u$, we obtain%
\[
\left(  \text{the }\left(  u,j\right)  \text{-th entry of the matrix
}A\right)  =a_{u,j}.
\]
Hence,%
\begin{align*}
a_{u,j}  &  =\left(  \text{the }\left(  u,j\right)  \text{-th entry of the
matrix }A\right) \\
&  =\left(  \text{the }j\text{-th entry of }\underbrace{\text{the }u\text{-th
row of }A}_{=\underbrace{\left(  0,0,\ldots,0\right)  }_{n\text{ zeroes}}%
}\right) \\
&  =\left(  \text{the }j\text{-th entry of }\underbrace{\left(  0,0,\ldots
,0\right)  }_{n\text{ zeroes}}\right)  =0.
\end{align*}
This proves (\ref{sol.ps4.6.c.1}).}. Now, every $\sigma\in S_{n}$ satisfies%
\[
\prod_{i=1}^{n}a_{i,\sigma\left(  i\right)  }=0
\]
\footnote{\textit{Proof.} Let $\sigma\in S_{n}$. Then, $a_{u,\sigma\left(
u\right)  }=0$ (by (\ref{sol.ps4.6.c.1})). But $a_{u,\sigma\left(  u\right)
}$ is a factor of the product $\prod_{i=1}^{n}a_{i,\sigma\left(  i\right)  }$
(namely, the factor for $i=u$). Hence, one factor of the product $\prod
_{i=1}^{n}a_{i,\sigma\left(  i\right)  }$ equals $0$ (since $a_{u,\sigma
\left(  u\right)  }=0$). Therefore, the whole product $\prod_{i=1}%
^{n}a_{i,\sigma\left(  i\right)  }$ equals $0$ (because if one factor of a
product equals $0$, then the whole product equals $0$). Qed.}. Now,
(\ref{eq.det.eq.2}) shows that%
\[
\det A=\sum_{\sigma\in S_{n}}\left(  -1\right)  ^{\sigma}\underbrace{\prod
_{i=1}^{n}a_{i,\sigma\left(  i\right)  }}_{=0}=\sum_{\sigma\in S_{n}}\left(
-1\right)  ^{\sigma}0=0.
\]
This solves Exercise \ref{exe.ps4.6} \textbf{(c)}.
\end{verlong}

\textbf{(d)} Assume that a column of $A$ consists of zeroes. Thus, a row of
$A^{T}$ consists of zeroes (because the columns of $A$ correspond to the rows
of $A^{T}$\ \ \ \ \footnote{See Remark \ref{rmk.colsA=rowsATT.rmk} for the
meaning of \textquotedblleft correspond\textquotedblright\ we are using
here.}). Hence, Exercise \ref{exe.ps4.6} \textbf{(c)} (applied to $A^{T}$
instead of $A$) yields $\det\left(  A^{T}\right)  =0$.

But Exercise \ref{exe.ps4.4} yields $\det\left(  A^{T}\right)  =\det A$.
Hence, $\det A=\det\left(  A^{T}\right)  =0$. This solves Exercise
\ref{exe.ps4.6} \textbf{(d)}.

\textbf{(e)} Assume that $A$ has two equal rows. In other words, some two
distinct rows of $A$ are equal (where \textquotedblleft
distinct\textquotedblright\ means that these rows are in different positions,
not that they are distinct as row vectors). In other words, there exist two
distinct elements $u$ and $v$ of $\left\{  1,2,\ldots,n\right\}  $ such that%
\begin{equation}
\left(  \text{the }u\text{-th row of }A\right)  =\left(  \text{the }v\text{-th
row of }A\right)  . \label{sol.ps4.6.e.1}%
\end{equation}
Consider these $u$ and $v$.

\begin{vershort}
The equality (\ref{sol.ps4.6.e.1}) shows that%
\begin{equation}
a_{u,j}=a_{v,j}\ \ \ \ \ \ \ \ \ \ \text{for every }j\in\left\{
1,2,\ldots,n\right\}  \label{sol.ps4.6.e.short.1a}%
\end{equation}
(because $a_{u,j}$ is the $j$-th entry of the $u$-th row of $A$, while
$a_{v,j}$ is the $j$-th entry of the $v$-th row of $A$).

Define a map $\kappa:\left[  n\right]  \rightarrow\left[  n\right]  $ by
\[
\left(  \kappa\left(  i\right)  =
\begin{cases}
i, & \text{if }i\neq u;\\
v, & \text{if }i=u
\end{cases}
\ \ \ \ \ \ \ \ \ \ \text{for every }i\in\left[  n\right]  \right)  .
\]
The definition of $\kappa$ shows that $\kappa\left(  v\right)  =v$ (since
$v\neq u$) but also $\kappa\left(  u\right)  =v$. Thus, $\kappa\left(
u\right)  =v=\kappa\left(  v\right)  $, in spite of $u\neq v$. Therefore, the
map $\kappa$ is not injective, and thus not bijective; in particular, $\kappa$
is not a permutation. Thus, $\kappa\notin S_{n}$. But every $i\in\left\{
1,2,\ldots,n\right\}  $ and $j\in\left\{  1,2,\ldots,n\right\}  $ satisfy%
\begin{equation}
a_{\kappa\left(  i\right)  ,j}=a_{i,j} \label{sol.ps4.6.e.short.3}%
\end{equation}
\footnote{\textit{Proof of (\ref{sol.ps4.6.e.short.3}):} Let $i\in\left\{
1,2,\ldots,n\right\}  $ and $j\in\left\{  1,2,\ldots,n\right\}  $. We must
prove (\ref{sol.ps4.6.e.short.3}). We must be in one of the following two
cases:
\par
\textit{Case 1:} We have $i\neq u$.
\par
\textit{Case 2:} We have $i=u$.
\par
Let us first consider Case 1. In this case, we have $i\neq u$. Thus, the
definition of $\kappa$ yields $\kappa\left(  i\right)  =
\begin{cases}
i, & \text{if }i\neq u;\\
v, & \text{if }i=u
\end{cases}
=i$ (since $i\neq u$), so that $a_{\kappa\left(  i\right)  ,j}=a_{i,j}$. Thus,
(\ref{sol.ps4.6.e.short.3}) is proven in Case 1.
\par
Let us now consider Case 2. In this case, we have $i=u$. Thus, $\kappa\left(
i\right)  =\kappa\left(  u\right)  =v$ (since $i=u$). Hence,%
\begin{align*}
a_{\kappa\left(  i\right)  ,j}  &  =a_{v,j}=a_{u,j}\ \ \ \ \ \ \ \ \ \ \left(
\text{by (\ref{sol.ps4.6.e.short.1a})}\right) \\
&  =a_{i,j}\ \ \ \ \ \ \ \ \ \ \left(  \text{since }u=i\right)  .
\end{align*}
Thus, (\ref{sol.ps4.6.e.short.3}) is proven in Case 2.
\par
We have now proven (\ref{sol.ps4.6.e.short.3}) in both Cases 1 and 2. Thus,
(\ref{sol.ps4.6.e.short.3}) always holds, qed.}. Thus, $\left(
\underbrace{a_{\kappa\left(  i\right)  ,j}}_{=a_{i,j}}\right)  _{1\leq i\leq
n,\ 1\leq j\leq n}=\left(  a_{i,j}\right)  _{1\leq i\leq n,\ 1\leq j\leq n}%
=A$. Therefore, we can apply Lemma \ref{lem.det.sigma} \textbf{(b)} to $A$,
$a_{i,j}$ and $A$ instead of $B$, $b_{i,j}$ and $B_{\kappa}$. We thus obtain
$\det A=0$. Exercise \ref{exe.ps4.6} \textbf{(e)} is thus solved.
\end{vershort}

\begin{verlong}
Define a map $\kappa:\left[  n\right]  \rightarrow\left[  n\right]  $ by
\[
\left(  \kappa\left(  i\right)  =
\begin{cases}
i, & \text{if }i\neq u;\\
v, & \text{if }i=u
\end{cases}
\ \ \ \ \ \ \ \ \ \ \text{for every }i\in\left[  n\right]  \right)  .
\]
Then, $\kappa\notin S_{n}$\ \ \ \ \footnote{\textit{Proof.} Assume the
contrary. Then, $\kappa\in S_{n}$. Thus, $\kappa$ is a permutation. Hence, the
map $\kappa$ is injective. But the definition of $\kappa$ yields
$\kappa\left(  u\right)  =
\begin{cases}
u, & \text{if }u\neq u;\\
v, & \text{if }u=u
\end{cases}
=v$ (since $u=u$). Also, $u$ and $v$ are distinct, so that we have $v\neq u$.
Hence, $\kappa\left(  v\right)  \neq\kappa\left(  u\right)  $ (since $\kappa$
is injective). But the definition of $\kappa$ yields $\kappa\left(  v\right)
=
\begin{cases}
v, & \text{if }v\neq u;\\
v, & \text{if }v=u
\end{cases}
=v=\kappa\left(  u\right)  $. This contradicts $\kappa\left(  v\right)
\neq\kappa\left(  u\right)  $. This contradiction proves that our assumption
was wrong, qed.}. But every $i\in\left\{  1,2,\ldots,n\right\}  $ and
$j\in\left\{  1,2,\ldots,n\right\}  $ satisfy%
\begin{equation}
a_{\kappa\left(  i\right)  ,j}=a_{i,j} \label{sol.ps4.6.e.3}%
\end{equation}
\footnote{\textit{Proof of (\ref{sol.ps4.6.e.3}):} Let $i\in\left\{
1,2,\ldots,n\right\}  $ and $j\in\left\{  1,2,\ldots,n\right\}  $. We must
prove (\ref{sol.ps4.6.e.3}). We must be in one of the following two cases:
\par
\textit{Case 1:} We have $i\neq u$.
\par
\textit{Case 2:} We have $i=u$.
\par
Let us first consider Case 1. In this case, we have $i\neq u$. Thus, the
definition of $\kappa$ yields $\kappa\left(  i\right)  =
\begin{cases}
i, & \text{if }i\neq u;\\
v, & \text{if }i=u
\end{cases}
=i$ (since $i\neq u$), so that $a_{\kappa\left(  i\right)  ,j}=a_{i,j}$. Thus,
(\ref{sol.ps4.6.e.3}) is proven in Case 1.
\par
Let us now consider Case 2. In this case, we have $i=u$. Thus, the definition
of $\kappa$ yields $\kappa\left(  i\right)  =
\begin{cases}
i, & \text{if }i\neq u;\\
v, & \text{if }i=u
\end{cases}
=v$ (since $i=u$). Hence,%
\begin{align*}
a_{\kappa\left(  i\right)  ,j}  &  =a_{v,j}=\left(  \text{the }\left(
v,j\right)  \text{-th entry of the matrix }A\right) \\
&  \ \ \ \ \ \ \ \ \ \ \left(
\begin{array}
[c]{c}%
\text{since }\left(  \text{the }\left(  v,j\right)  \text{-th entry of the
matrix }A\right)  =a_{v,j}\\
\text{(by (\ref{sol.ps4.6.aij}), applied to }v\text{ instead of }i\text{)}%
\end{array}
\right) \\
&  =\left(  \text{the }j\text{-th entry of }\underbrace{\text{the }v\text{-th
row of }A}_{\substack{=\left(  \text{the }u\text{-th row of }A\right)
\\\text{(by (\ref{sol.ps4.6.e.1}))}}}\right)  =\left(  \text{the }j\text{-th
entry of the }u\text{-th row of }A\right) \\
&  =\left(  \text{the }\left(  u,j\right)  \text{-th entry of the matrix
}A\right)  =a_{u,j}\ \ \ \ \ \ \ \ \ \ \left(  \text{by (\ref{sol.ps4.6.aij}),
applied to }u\text{ instead of }i\right) \\
&  =a_{i,j}\ \ \ \ \ \ \ \ \ \ \left(  \text{since }u=i\right)  .
\end{align*}
Thus, (\ref{sol.ps4.6.e.3}) is proven in Case 2.
\par
We have now proven (\ref{sol.ps4.6.e.3}) in both Cases 1 and 2. Thus,
(\ref{sol.ps4.6.e.3}) always holds, qed.}. Thus, $\left(
\underbrace{a_{\kappa\left(  i\right)  ,j}}_{=a_{i,j}}\right)  _{1\leq i\leq
n,\ 1\leq j\leq n}=\left(  a_{i,j}\right)  _{1\leq i\leq n,\ 1\leq j\leq n}%
=A$. Therefore, we can apply Lemma \ref{lem.det.sigma} \textbf{(b)} to $A$,
$a_{i,j}$ and $A$ instead of $B$, $b_{i,j}$ and $B_{\kappa}$. We thus obtain
$\det A=0$. Exercise \ref{exe.ps4.6} \textbf{(e)} is thus solved.
\end{verlong}

\textbf{(f)} Assume that $A$ has two equal columns. Thus, $A^{T}$ has two
equal rows (because the columns of $A$ correspond to the rows of $A^{T}%
$\ \ \ \ \footnote{See Remark \ref{rmk.colsA=rowsATT.rmk} for the meaning of
\textquotedblleft correspond\textquotedblright\ we are using here.}). Hence,
Exercise \ref{exe.ps4.6} \textbf{(e)} (applied to $A^{T}$ instead of $A$)
yields $\det\left(  A^{T}\right)  =0$.

But Exercise \ref{exe.ps4.4} yields $\det\left(  A^{T}\right)  =\det A$.
Hence, $\det A=\det\left(  A^{T}\right)  =0$. This solves Exercise
\ref{exe.ps4.6} \textbf{(f)}.

\textbf{(g)} Let $B$ be the $n\times n$-matrix obtained from $A$ by
multiplying the $k$-th row by $\lambda$. Thus,
\begin{equation}
\left(  \text{the }k\text{-th row of }B\right)  =\lambda\left(  \text{the
}k\text{-th row of }A\right)  , \label{sol.ps4.6.g.krow}%
\end{equation}
whereas%
\begin{align}
&  \left(  \left(  \text{the }u\text{-th row of }B\right)  =\left(  \text{the
}u\text{-th row of }A\right)  \right. \label{sol.ps4.6.g.urow}\\
&  \ \ \ \ \ \ \ \ \ \ \left.  \text{for all }u\in\left\{  1,2,\ldots
,n\right\}  \text{ satisfying }u\neq k\right)  .\nonumber
\end{align}

We write the matrix $B$ in the form $B=\left(  b_{i,j}\right)  _{1\leq i\leq
n,\ 1\leq j\leq n}$. Thus, every $\left(  i,j\right)  \in\left\{
1,2,\ldots,n\right\}  ^{2}$ satisfies%
\begin{equation}
\left(  \text{the }\left(  i,j\right)  \text{-th entry of the matrix
}B\right)  =b_{i,j}. \label{sol.ps4.6.g.bij}%
\end{equation}

\begin{vershort}
For every $u\in\left\{  1,2,\ldots,n\right\}  $ and $v\in\left\{
1,2,\ldots,n\right\}  $ satisfying $u\neq k$, we have%
\begin{equation}
b_{u,v}=a_{u,v} \label{sol.ps4.6.g.short.buv}%
\end{equation}
\footnote{\textit{Proof of (\ref{sol.ps4.6.g.short.buv}):} Let $u\in\left\{
1,2,\ldots,n\right\}  $ and $v\in\left\{  1,2,\ldots,n\right\}  $ be such that
$u\neq k$. Then, (\ref{sol.ps4.6.g.urow}) shows that the $v$-th entry of the
$u$-th row of $B$ equals the $v$-th entry of the $u$-th row of $A$. Since the
former entry is $b_{u,v}$, while the latter entry equals $a_{u,v}$, this
rewrites as $b_{u,v}=a_{u,v}$. This proves (\ref{sol.ps4.6.g.short.buv}).}.
For every $v\in\left\{  1,2,\ldots,n\right\}  $, we have
\begin{equation}
b_{k,v}=\lambda a_{k,v} \label{sol.ps4.6.g.short.bkv}%
\end{equation}
\footnote{\textit{Proof of (\ref{sol.ps4.6.g.short.bkv}):} Let $v\in\left\{
1,2,\ldots,n\right\}  $. Then, (\ref{sol.ps4.6.g.krow}) shows that the $v$-th
entry of the $k$-th row of $B$ equals $\lambda$ times the $v$-th entry of the
$k$-th row of $A$. Since the former entry is $b_{k,v}$, while the latter entry
equals $a_{k,v}$, this rewrites as $b_{k,v}=\lambda a_{k,v}$. This proves
(\ref{sol.ps4.6.g.short.bkv}).}. Now, it is easy to see that every $\sigma\in
S_{n}$ satisfies%
\begin{equation}
\prod_{i=1}^{n}b_{i,\sigma\left(  i\right)  }=\lambda\prod_{i=1}%
^{n}a_{i,\sigma\left(  i\right)  } \label{sol.ps4.6.g.short.sigma}%
\end{equation}
\footnote{\textit{Proof of (\ref{sol.ps4.6.g.short.sigma}):} Let $\sigma\in
S_{n}$. Taking the $k$-th factor out of the product $\prod_{i=1}%
^{n}a_{i,\sigma\left(  i\right)  }$, we obtain
\[
\prod_{i=1}^{n}a_{i,\sigma\left(  i\right)  }=a_{k,\sigma\left(  k\right)
}\cdot\prod_{\substack{i\in\left\{  1,2,\ldots,n\right\}  ;\\i\neq
k}}a_{i,\sigma\left(  i\right)  }.
\]
Multiplying both sides of this equality by $\lambda$, we obtain%
\[
\lambda\prod_{i=1}^{n}a_{i,\sigma\left(  i\right)  }=\lambda a_{k,\sigma
\left(  k\right)  }\cdot\prod_{\substack{i\in\left\{  1,2,\ldots,n\right\}
;\\i\neq k}}a_{i,\sigma\left(  i\right)  }.
\]
Compared with%
\begin{align*}
\prod_{i=1}^{n}b_{i,\sigma\left(  i\right)  }  &  =\underbrace{b_{k,\sigma
\left(  k\right)  }}_{\substack{=\lambda a_{k,\sigma\left(  k\right)
}\\\text{(by (\ref{sol.ps4.6.g.short.bkv}), applied to}\\v=\sigma\left(
k\right)  \text{)}}}\cdot\prod_{\substack{i\in\left\{  1,2,\ldots,n\right\}
;\\i\neq k}}\underbrace{b_{i,\sigma\left(  i\right)  }}%
_{\substack{=a_{i,\sigma\left(  i\right)  }\\\text{(by
(\ref{sol.ps4.6.g.short.buv}), applied to}\\u=i\text{ and }v=\sigma\left(
i\right)  \text{)}}}\ \ \ \ \ \ \ \ \ \ \left(  \text{since }k\in\left\{
1,2,\ldots,n\right\}  \right) \\
&  =\lambda a_{k,\sigma\left(  k\right)  }\cdot\prod_{\substack{i\in\left\{
1,2,\ldots,n\right\}  ;\\i\neq k}}a_{i,\sigma\left(  i\right)  },
\end{align*}
this yields $\prod_{i=1}^{n}b_{i,\sigma\left(  i\right)  }=\lambda\prod
_{i=1}^{n}a_{i,\sigma\left(  i\right)  }$. This proves
(\ref{sol.ps4.6.g.short.sigma}).}. Now, recall that $B=\left(  b_{i,j}\right)
_{1\leq i\leq n,\ 1\leq j\leq n}$. Hence, (\ref{eq.det.eq.2}) (applied to $B$
and $b_{i,j}$ instead of $A$ and $a_{i,j}$) yields%
\begin{align*}
\det B  &  =\sum_{\sigma\in S_{n}}\left(  -1\right)  ^{\sigma}%
\underbrace{\prod_{i=1}^{n}b_{i,\sigma\left(  i\right)  }}_{\substack{=\lambda
\prod_{i=1}^{n}a_{i,\sigma\left(  i\right)  }\\\text{(by
(\ref{sol.ps4.6.g.short.sigma}))}}}=\sum_{\sigma\in S_{n}}\left(  -1\right)
^{\sigma}\lambda\prod_{i=1}^{n}a_{i,\sigma\left(  i\right)  }\\
&  =\lambda\underbrace{\sum_{\sigma\in S_{n}}\left(  -1\right)  ^{\sigma}%
\prod_{i=1}^{n}a_{i,\sigma\left(  i\right)  }}_{\substack{=\det A\\\text{(by
(\ref{eq.det.eq.2}))}}}=\lambda\det A.
\end{align*}
Thus, Exercise \ref{exe.ps4.6} \textbf{(g)} is solved.
\end{vershort}

\begin{verlong}
For every $u\in\left\{  1,2,\ldots,n\right\}  $ and $v\in\left\{
1,2,\ldots,n\right\}  $ satisfying $u\neq k$, we have%
\begin{equation}
b_{u,v}=a_{u,v} \label{sol.ps4.6.g.buv}%
\end{equation}
\footnote{\textit{Proof of (\ref{sol.ps4.6.g.buv}):} Let $u\in\left\{
1,2,\ldots,n\right\}  $ and $v\in\left\{  1,2,\ldots,n\right\}  $ be such that
$u\neq k$. Then, (\ref{sol.ps4.6.g.bij}) (applied to $i=u$ and $j=v$) yields%
\[
\left(  \text{the }\left(  u,v\right)  \text{-th entry of the matrix
}B\right)  =b_{u,v},
\]
so that%
\begin{align*}
b_{u,v}  &  =\left(  \text{the }\left(  u,v\right)  \text{-th entry of the
matrix }B\right) \\
&  =\left(  \text{the }v\text{-th entry of }\underbrace{\text{the }u\text{-th
row of }B}_{\substack{=\left(  \text{the }u\text{-th row of }A\right)
\\\text{(by (\ref{sol.ps4.6.g.urow}))}}}\right) \\
&  =\left(  \text{the }v\text{-th entry of the }u\text{-th row of }A\right) \\
&  =\left(  \text{the }\left(  u,v\right)  \text{-th entry of the matrix
}A\right)  =a_{u,v}%
\end{align*}
(by (\ref{sol.ps4.6.aij}), applied to $i=u$ and $j=v$), qed.}. For every
$v\in\left\{  1,2,\ldots,n\right\}  $, we have
\begin{equation}
b_{k,v}=\lambda a_{k,v} \label{sol.ps4.6.g.bkv}%
\end{equation}
\footnote{\textit{Proof of (\ref{sol.ps4.6.g.bkv}):} Let $v\in\left\{
1,2,\ldots,n\right\}  $. Then, (\ref{sol.ps4.6.g.bij}) (applied to $i=k$ and
$j=v$) yields%
\[
\left(  \text{the }\left(  k,v\right)  \text{-th entry of the matrix
}B\right)  =b_{k,v},
\]
so that%
\begin{align*}
b_{k,v}  &  =\left(  \text{the }\left(  k,v\right)  \text{-th entry of the
matrix }B\right) \\
&  =\left(  \text{the }v\text{-th entry of }\underbrace{\text{the }k\text{-th
row of }B}_{=\lambda\left(  \text{the }k\text{-th row of }A\right)  }\right)
\\
&  =\left(  \text{the }v\text{-th entry of }\lambda\left(  \text{the
}k\text{-th row of }A\right)  \right) \\
&  =\lambda\underbrace{\left(  \text{the }v\text{-th entry of the }k\text{-th
row of }A\right)  }_{=\left(  \text{the }\left(  k,v\right)  \text{-th entry
of the matrix }A\right)  }\\
&  =\lambda\underbrace{\left(  \text{the }\left(  k,v\right)  \text{-th entry
of the matrix }A\right)  }_{\substack{=a_{k,v}\\\text{(by (\ref{sol.ps4.6.aij}%
), applied to }i=k\text{ and }j=v\text{)}}}=\lambda a_{k,v},
\end{align*}
qed.}. Now, it is easy to see that every $\sigma\in S_{n}$ satisfies%
\begin{equation}
\prod_{i=1}^{n}b_{i,\sigma\left(  i\right)  }=\lambda\prod_{i=1}%
^{n}a_{i,\sigma\left(  i\right)  } \label{sol.ps4.6.g.sigma}%
\end{equation}
\footnote{\textit{Proof of (\ref{sol.ps4.6.g.sigma}):} Let $\sigma\in S_{n}$.
We have%
\[
\underbrace{\prod_{i=1}^{n}}_{=\prod_{i\in\left\{  1,2,\ldots,n\right\}  }%
}a_{i,\sigma\left(  i\right)  }=\prod_{i\in\left\{  1,2,\ldots,n\right\}
}a_{i,\sigma\left(  i\right)  }=a_{k,\sigma\left(  k\right)  }\cdot
\prod_{\substack{i\in\left\{  1,2,\ldots,n\right\}  ;\\i\neq k}}a_{i,\sigma
\left(  i\right)  }%
\]
(since $k\in\left\{  1,2,\ldots,n\right\}  $). Multiplying both sides of this
equality by $\lambda$, we obtain%
\[
\lambda\prod_{i=1}^{n}a_{i,\sigma\left(  i\right)  }=\lambda a_{k,\sigma
\left(  k\right)  }\cdot\prod_{\substack{i\in\left\{  1,2,\ldots,n\right\}
;\\i\neq k}}a_{i,\sigma\left(  i\right)  }.
\]
Compared with%
\begin{align*}
\underbrace{\prod_{i=1}^{n}}_{=\prod_{i\in\left\{  1,2,\ldots,n\right\}  }%
}b_{i,\sigma\left(  i\right)  }  &  =\prod_{i\in\left\{  1,2,\ldots,n\right\}
}b_{i,\sigma\left(  i\right)  }\\
&  =\underbrace{b_{k,\sigma\left(  k\right)  }}_{\substack{=\lambda
a_{k,\sigma\left(  k\right)  }\\\text{(by (\ref{sol.ps4.6.g.bkv}), applied
to}\\v=\sigma\left(  k\right)  \text{)}}}\cdot\prod_{\substack{i\in\left\{
1,2,\ldots,n\right\}  ;\\i\neq k}}\underbrace{b_{i,\sigma\left(  i\right)  }%
}_{\substack{=a_{i,\sigma\left(  i\right)  }\\\text{(by (\ref{sol.ps4.6.g.buv}%
), applied to}\\u=i\text{ and }v=\sigma\left(  i\right)  \text{)}%
}}\ \ \ \ \ \ \ \ \ \ \left(  \text{since }k\in\left\{  1,2,\ldots,n\right\}
\right) \\
&  =\lambda a_{k,\sigma\left(  k\right)  }\cdot\prod_{\substack{i\in\left\{
1,2,\ldots,n\right\}  ;\\i\neq k}}a_{i,\sigma\left(  i\right)  },
\end{align*}
this yields $\prod_{i=1}^{n}b_{i,\sigma\left(  i\right)  }=\lambda\prod
_{i=1}^{n}a_{i,\sigma\left(  i\right)  }$. This proves
(\ref{sol.ps4.6.g.sigma}).}. Now, recall that $B=\left(  b_{i,j}\right)
_{1\leq i\leq n,\ 1\leq j\leq n}$. Hence, (\ref{eq.det.eq.2}) (applied to $B$
and $b_{i,j}$ instead of $A$ and $a_{i,j}$) yields%
\begin{align*}
\det B  &  =\sum_{\sigma\in S_{n}}\left(  -1\right)  ^{\sigma}%
\underbrace{\prod_{i=1}^{n}b_{i,\sigma\left(  i\right)  }}_{\substack{=\lambda
\prod_{i=1}^{n}a_{i,\sigma\left(  i\right)  }\\\text{(by
(\ref{sol.ps4.6.g.sigma}))}}}=\sum_{\sigma\in S_{n}}\left(  -1\right)
^{\sigma}\lambda\prod_{i=1}^{n}a_{i,\sigma\left(  i\right)  }\\
&  =\lambda\underbrace{\sum_{\sigma\in S_{n}}\left(  -1\right)  ^{\sigma}%
\prod_{i=1}^{n}a_{i,\sigma\left(  i\right)  }}_{\substack{=\det A\\\text{(by
(\ref{eq.det.eq.2}))}}}=\lambda\det A.
\end{align*}
Thus, Exercise \ref{exe.ps4.6} \textbf{(g)} is solved.
\end{verlong}

\textbf{(h)} Let $B$ be the $n\times n$-matrix obtained from $A$ by
multiplying the $k$-th column by $\lambda$. Thus, $B^{T}$ is the $n\times
n$-matrix obtained from $A^{T}$ by multiplying the $k$-th row by $\lambda$
(because the columns of $A$ correspond to the rows of $A^{T}$%
\ \ \ \ \footnote{See Remark \ref{rmk.colsA=rowsATT.rmk} for the meaning of
\textquotedblleft correspond\textquotedblright\ we are using here.}). Hence,
Exercise \ref{exe.ps4.6} \textbf{(g)} (applied to $A^{T}$ and $B^{T}$ instead
of $A$ and $B$) yields $\det\left(  B^{T}\right)  =\lambda\det\left(
A^{T}\right)  $.

But Exercise \ref{exe.ps4.4} yields $\det\left(  A^{T}\right)  =\det A$. Also,
Exercise \ref{exe.ps4.4} (applied to $B$ instead of $A$) yields $\det\left(
B^{T}\right)  =\det B$. But recall that $\det\left(  B^{T}\right)
=\lambda\det\left(  A^{T}\right)  $. This rewrites as $\det B=\lambda\det A$
(since $\det\left(  B^{T}\right)  =\det B$ and $\det\left(  A^{T}\right)
=\det A$). This solves Exercise \ref{exe.ps4.6} \textbf{(h)}.

\textbf{(i)} We know that the rows of the matrix $A^{\prime}$ equal the
corresponding rows of $A$ except (perhaps) the $k$-th row. In other words,%
\begin{align}
&  \left(  \left(  \text{the }u\text{-th row of }A^{\prime}\right)  =\left(
\text{the }u\text{-th row of }A\right)  \right. \label{sol.ps4.6.i.urowA'}\\
&  \ \ \ \ \ \ \ \ \ \ \left.  \text{for all }u\in\left\{  1,2,\ldots
,n\right\}  \text{ satisfying }u\neq k\right)  .\nonumber
\end{align}

We know that $B$ is the $n\times n$-matrix obtained from $A$ by adding the
$k$-th row of $A^{\prime}$ to the $k$-th row of $A$. Hence,%
\begin{equation}
\left(  \text{the }k\text{-th row of }B\right)  =\left(  \text{the }k\text{-th
row of }A\right)  +\left(  \text{the }k\text{-th row of }A^{\prime}\right)
\label{sol.ps4.6.i.krowB}%
\end{equation}
and%
\begin{align}
&  \left(  \left(  \text{the }u\text{-th row of }B\right)  =\left(  \text{the
}u\text{-th row of }A\right)  \right. \label{sol.ps4.6.i.urowB}\\
&  \ \ \ \ \ \ \ \ \ \ \left.  \text{for all }u\in\left\{  1,2,\ldots
,n\right\}  \text{ satisfying }u\neq k\right)  .\nonumber
\end{align}

\begin{vershort}
We write the matrix $B$ in the form $B=\left(  b_{i,j}\right)  _{1\leq i\leq
n,\ 1\leq j\leq n}$, and we write the matrix $A^{\prime}$ in the form
$A^{\prime}=\left(  a_{i,j}^{\prime}\right)  _{1\leq i\leq n,\ 1\leq j\leq n}%
$. Then, for every $u\in\left\{  1,2,\ldots,n\right\}  $ and $v\in\left\{
1,2,\ldots,n\right\}  $ satisfying $u\neq k$, we have%
\begin{equation}
a_{u,v}^{\prime}=a_{u,v} \label{sol.ps4.6.i.short.a'auv}%
\end{equation}
\footnote{This follows from (\ref{sol.ps4.6.i.urowA'}) in a similar way as
(\ref{sol.ps4.6.g.short.buv}) follows from (\ref{sol.ps4.6.g.urow}).} and%
\begin{equation}
b_{u,v}=a_{u,v} \label{sol.ps4.6.i.short.bauv}%
\end{equation}
\footnote{This follows from (\ref{sol.ps4.6.i.urowB}) in a similar way as
(\ref{sol.ps4.6.g.short.buv}) follows from (\ref{sol.ps4.6.g.urow}).}. Also,
for every $v\in\left\{  1,2,\ldots,n\right\}  $, we have
\begin{equation}
b_{k,v}=a_{k,v}+a_{k,v}^{\prime} \label{sol.ps4.6.i.short.bkv}%
\end{equation}
\footnote{\textit{Proof of (\ref{sol.ps4.6.i.short.bkv}):} Let $v\in\left\{
1,2,\ldots,n\right\}  $. Then, (\ref{sol.ps4.6.i.krowB}) shows that
\begin{align*}
&  \left(  \text{the }v\text{-th entry of the }k\text{-th row of }B\right) \\
&  =\left(  \text{the }v\text{-th entry of the }k\text{-th row of }A\right)
+\left(  \text{the }v\text{-th entry of the }k\text{-th row of }A^{\prime
}\right)  .
\end{align*}
But the three entries appearing in this equality are $b_{k,v}$, $a_{k,v}$ and
$a_{k,v}^{\prime}$ (in this order). Thus, this equality rewrites as
$b_{k,v}=a_{k,v}+a_{k,v}^{\prime}$. This proves (\ref{sol.ps4.6.i.short.bkv}%
).}. Now, it is easy to see that every $\sigma\in S_{n}$ satisfies%
\begin{equation}
\prod_{i=1}^{n}b_{i,\sigma\left(  i\right)  }=\prod_{i=1}^{n}a_{i,\sigma
\left(  i\right)  }+\prod_{i=1}^{n}a_{i,\sigma\left(  i\right)  }^{\prime}
\label{sol.ps4.6.i.short.sigma}%
\end{equation}
\footnote{\textit{Proof of (\ref{sol.ps4.6.i.short.sigma}):} Let $\sigma\in
S_{n}$. Pulling the $k$-th factor out of the product $\prod_{i=1}%
^{n}a_{i,\sigma\left(  i\right)  }$, we obtain%
\begin{equation}
\prod_{i=1}^{n}a_{i,\sigma\left(  i\right)  }=a_{k,\sigma\left(  k\right)
}\cdot\prod_{\substack{i\in\left\{  1,2,\ldots,n\right\}  ;\\i\neq
k}}a_{i,\sigma\left(  i\right)  } \label{sol.ps4.6.i.short.sigma.pf.1}%
\end{equation}
(since $k\in\left\{  1,2,\ldots,n\right\}  $). Similarly,%
\begin{equation}
\prod_{i=1}^{n}a_{i,\sigma\left(  i\right)  }^{\prime}=a_{k,\sigma\left(
k\right)  }^{\prime}\cdot\prod_{\substack{i\in\left\{  1,2,\ldots,n\right\}
;\\i\neq k}}\underbrace{a_{i,\sigma\left(  i\right)  }^{\prime}}%
_{\substack{=a_{i,\sigma\left(  i\right)  }\\\text{(by
(\ref{sol.ps4.6.i.short.a'auv}), applied to}\\u=i\text{ and }v=\sigma\left(
i\right)  \text{)}}}=a_{k,\sigma\left(  k\right)  }^{\prime}\cdot
\prod_{\substack{i\in\left\{  1,2,\ldots,n\right\}  ;\\i\neq k}}a_{i,\sigma
\left(  i\right)  } \label{sol.ps4.6.i.short.sigma.pf.2}%
\end{equation}
and%
\begin{align*}
\prod_{i=1}^{n}b_{i,\sigma\left(  i\right)  }  &  =\underbrace{b_{k,\sigma
\left(  k\right)  }}_{\substack{=a_{k,\sigma\left(  k\right)  }+a_{k,\sigma
\left(  k\right)  }^{\prime}\\\text{(by (\ref{sol.ps4.6.i.short.bkv}), applied
to}\\v=\sigma\left(  k\right)  \text{)}}}\cdot\prod_{\substack{i\in\left\{
1,2,\ldots,n\right\}  ;\\i\neq k}}\underbrace{b_{i,\sigma\left(  i\right)  }%
}_{\substack{=a_{i,\sigma\left(  i\right)  }\\\text{(by
(\ref{sol.ps4.6.i.short.bauv}), applied to}\\u=i\text{ and }v=\sigma\left(
i\right)  \text{)}}}\\
&  =\left(  a_{k,\sigma\left(  k\right)  }+a_{k,\sigma\left(  k\right)
}^{\prime}\right)  \cdot\prod_{\substack{i\in\left\{  1,2,\ldots,n\right\}
;\\i\neq k}}a_{i,\sigma\left(  i\right)  }\\
&  =\underbrace{a_{k,\sigma\left(  k\right)  }\cdot\prod_{\substack{i\in
\left\{  1,2,\ldots,n\right\}  ;\\i\neq k}}a_{i,\sigma\left(  i\right)  }%
}_{\substack{=\prod_{i=1}^{n}a_{i,\sigma\left(  i\right)  }\\\text{(by
(\ref{sol.ps4.6.i.short.sigma.pf.1}))}}}+\underbrace{a_{k,\sigma\left(
k\right)  }^{\prime}\cdot\prod_{\substack{i\in\left\{  1,2,\ldots,n\right\}
;\\i\neq k}}a_{i,\sigma\left(  i\right)  }}_{\substack{=\prod_{i=1}%
^{n}a_{i,\sigma\left(  i\right)  }^{\prime}\\\text{(by
(\ref{sol.ps4.6.i.short.sigma.pf.2}))}}}\\
&  =\prod_{i=1}^{n}a_{i,\sigma\left(  i\right)  }+\prod_{i=1}^{n}%
a_{i,\sigma\left(  i\right)  }^{\prime}.
\end{align*}
This proves (\ref{sol.ps4.6.i.short.sigma}).}.
\end{vershort}

\begin{verlong}
We write the matrix $B$ in the form $B=\left(  b_{i,j}\right)  _{1\leq i\leq
n,\ 1\leq j\leq n}$. Thus, every $\left(  i,j\right)  \in\left\{
1,2,\ldots,n\right\}  ^{2}$ satisfies%
\begin{equation}
\left(  \text{the }\left(  i,j\right)  \text{-th entry of the matrix
}B\right)  =b_{i,j}. \label{sol.ps4.6.i.bij}%
\end{equation}
We write the matrix $A^{\prime}$ in the form $A^{\prime}=\left(
a_{i,j}^{\prime}\right)  _{1\leq i\leq n,\ 1\leq j\leq n}$. Thus, every
$\left(  i,j\right)  \in\left\{  1,2,\ldots,n\right\}  ^{2}$ satisfies%
\begin{equation}
\left(  \text{the }\left(  i,j\right)  \text{-th entry of the matrix
}A^{\prime}\right)  =a_{i,j}^{\prime}. \label{sol.ps4.6.i.a'ij}%
\end{equation}

For every $u\in\left\{  1,2,\ldots,n\right\}  $ and $v\in\left\{
1,2,\ldots,n\right\}  $ satisfying $u\neq k$, we have%
\begin{equation}
b_{u,v}=a_{u,v} \label{sol.ps4.6.i.buv}%
\end{equation}
\footnote{\textit{Proof of (\ref{sol.ps4.6.i.buv}):} Let $u\in\left\{
1,2,\ldots,n\right\}  $ and $v\in\left\{  1,2,\ldots,n\right\}  $ be such that
$u\neq k$. Then, (\ref{sol.ps4.6.i.bij}) (applied to $i=u$ and $j=v$) yields%
\[
\left(  \text{the }\left(  u,v\right)  \text{-th entry of the matrix
}B\right)  =b_{u,v},
\]
so that%
\begin{align*}
b_{u,v}  &  =\left(  \text{the }\left(  u,v\right)  \text{-th entry of the
matrix }B\right) \\
&  =\left(  \text{the }v\text{-th entry of }\underbrace{\text{the }u\text{-th
row of }B}_{\substack{=\left(  \text{the }u\text{-th row of }A\right)
\\\text{(by (\ref{sol.ps4.6.i.urowB}))}}}\right) \\
&  =\left(  \text{the }v\text{-th entry of the }u\text{-th row of }A\right) \\
&  =\left(  \text{the }\left(  u,v\right)  \text{-th entry of the matrix
}A\right)  =a_{u,v}%
\end{align*}
(by (\ref{sol.ps4.6.aij}), applied to $i=u$ and $j=v$), qed.} and%
\begin{equation}
b_{u,v}=a_{u,v}^{\prime} \label{sol.ps4.6.i.a'uv}%
\end{equation}
\footnote{\textit{Proof of (\ref{sol.ps4.6.i.a'uv}):} Let $u\in\left\{
1,2,\ldots,n\right\}  $ and $v\in\left\{  1,2,\ldots,n\right\}  $ be such that
$u\neq k$. Then, (\ref{sol.ps4.6.i.a'ij}) (applied to $i=u$ and $j=v$) yields%
\[
\left(  \text{the }\left(  u,v\right)  \text{-th entry of the matrix
}A^{\prime}\right)  =a_{u,v}^{\prime},
\]
so that%
\begin{align*}
a_{u,v}^{\prime}  &  =\left(  \text{the }\left(  u,v\right)  \text{-th entry
of the matrix }A^{\prime}\right) \\
&  =\left(  \text{the }v\text{-th entry of }\underbrace{\text{the }u\text{-th
row of }A^{\prime}}_{\substack{=\left(  \text{the }u\text{-th row of
}A\right)  \\\text{(by (\ref{sol.ps4.6.i.urowA'}))}}}\right) \\
&  =\left(  \text{the }v\text{-th entry of the }u\text{-th row of }A\right) \\
&  =\left(  \text{the }\left(  u,v\right)  \text{-th entry of the matrix
}A\right)  =a_{u,v}%
\end{align*}
(by (\ref{sol.ps4.6.aij}), applied to $i=u$ and $j=v$). Compared with
(\ref{sol.ps4.6.i.buv}), this yields $b_{u,v}=a_{u,v}^{\prime}$, qed.}. For
every $v\in\left\{  1,2,\ldots,n\right\}  $, we have
\begin{equation}
b_{k,v}=a_{k,v}+a_{k,v}^{\prime} \label{sol.ps4.6.i.bkv}%
\end{equation}
\footnote{\textit{Proof of (\ref{sol.ps4.6.i.bkv}):} Let $v\in\left\{
1,2,\ldots,n\right\}  $. Then, (\ref{sol.ps4.6.i.bij}) (applied to $i=k$ and
$j=v$) yields%
\[
\left(  \text{the }\left(  k,v\right)  \text{-th entry of the matrix
}B\right)  =b_{k,v},
\]
so that%
\begin{align*}
b_{k,v}  &  =\left(  \text{the }\left(  k,v\right)  \text{-th entry of the
matrix }B\right) \\
&  =\left(  \text{the }v\text{-th entry of }\underbrace{\text{the }k\text{-th
row of }B}_{=\left(  \text{the }k\text{-th row of }A\right)  +\left(
\text{the }k\text{-th row of }A^{\prime}\right)  }\right) \\
&  =\left(  \text{the }v\text{-th entry of }\left(  \text{the }k\text{-th row
of }A\right)  +\left(  \text{the }k\text{-th row of }A^{\prime}\right)
\right) \\
&  =\underbrace{\left(  \text{the }v\text{-th entry of the }k\text{-th row of
}A\right)  }_{=\left(  \text{the }\left(  k,v\right)  \text{-th entry of the
matrix }A\right)  }\\
&  \ \ \ \ \ \ \ \ \ \ +\underbrace{\left(  \text{the }v\text{-th entry of the
}k\text{-th row of }A^{\prime}\right)  }_{=\left(  \text{the }\left(
k,v\right)  \text{-th entry of the matrix }A^{\prime}\right)  }\\
&  =\underbrace{\left(  \text{the }\left(  k,v\right)  \text{-th entry of the
matrix }A\right)  }_{\substack{=a_{k,v}\\\text{(by (\ref{sol.ps4.6.aij}),
applied to }i=k\text{ and }j=v\text{)}}}+\underbrace{\left(  \text{the
}\left(  k,v\right)  \text{-th entry of the matrix }A^{\prime}\right)
}_{\substack{=a_{k,v}^{\prime}\\\text{(by (\ref{sol.ps4.6.i.a'ij}), applied to
}i=k\text{ and }j=v\text{)}}}\\
&  =a_{k,v}+a_{k,v}^{\prime},
\end{align*}
qed.}. Now, it is easy to see that every $\sigma\in S_{n}$ satisfies%
\begin{equation}
\prod_{i=1}^{n}b_{i,\sigma\left(  i\right)  }=\prod_{i=1}^{n}a_{i,\sigma
\left(  i\right)  }+\prod_{i=1}^{n}a_{i,\sigma\left(  i\right)  }^{\prime}
\label{sol.ps4.6.i.sigma}%
\end{equation}
\footnote{\textit{Proof of (\ref{sol.ps4.6.i.sigma}):} Let $\sigma\in S_{n}$.
We have%
\begin{equation}
\underbrace{\prod_{i=1}^{n}}_{=\prod_{i\in\left\{  1,2,\ldots,n\right\}  }%
}a_{i,\sigma\left(  i\right)  }=\prod_{i\in\left\{  1,2,\ldots,n\right\}
}a_{i,\sigma\left(  i\right)  }=a_{k,\sigma\left(  k\right)  }\cdot
\prod_{\substack{i\in\left\{  1,2,\ldots,n\right\}  ;\\i\neq k}}a_{i,\sigma
\left(  i\right)  } \label{sol.ps4.6.i.sigma.pf.1}%
\end{equation}
(since $k\in\left\{  1,2,\ldots,n\right\}  $) and%
\begin{equation}
\underbrace{\prod_{i=1}^{n}}_{=\prod_{i\in\left\{  1,2,\ldots,n\right\}  }%
}a_{i,\sigma\left(  i\right)  }^{\prime}=\prod_{i\in\left\{  1,2,\ldots
,n\right\}  }a_{i,\sigma\left(  i\right)  }^{\prime}=a_{k,\sigma\left(
k\right)  }^{\prime}\cdot\prod_{\substack{i\in\left\{  1,2,\ldots,n\right\}
;\\i\neq k}}a_{i,\sigma\left(  i\right)  }^{\prime}
\label{sol.ps4.6.i.sigma.pf.2}%
\end{equation}
(since $k\in\left\{  1,2,\ldots,n\right\}  $). Furthermore,%
\begin{align*}
\underbrace{\prod_{i=1}^{n}}_{=\prod_{i\in\left\{  1,2,\ldots,n\right\}  }%
}b_{i,\sigma\left(  i\right)  }  &  =\prod_{i\in\left\{  1,2,\ldots,n\right\}
}b_{i,\sigma\left(  i\right)  }=\underbrace{b_{k,\sigma\left(  k\right)  }%
}_{\substack{=a_{k,\sigma\left(  k\right)  }+a_{k,\sigma\left(  k\right)
}^{\prime}\\\text{(by (\ref{sol.ps4.6.i.bkv}), applied to}\\v=\sigma\left(
k\right)  \text{)}}}\cdot\prod_{\substack{i\in\left\{  1,2,\ldots,n\right\}
;\\i\neq k}}b_{i,\sigma\left(  i\right)  }\\
&  \ \ \ \ \ \ \ \ \ \ \left(  \text{since }k\in\left\{  1,2,\ldots,n\right\}
\right) \\
&  =\left(  a_{k,\sigma\left(  k\right)  }+a_{k,\sigma\left(  k\right)
}^{\prime}\right)  \cdot\prod_{\substack{i\in\left\{  1,2,\ldots,n\right\}
;\\i\neq k}}b_{i,\sigma\left(  i\right)  }\\
&  =a_{k,\sigma\left(  k\right)  }\cdot\prod_{\substack{i\in\left\{
1,2,\ldots,n\right\}  ;\\i\neq k}}\underbrace{b_{i,\sigma\left(  i\right)  }%
}_{\substack{=a_{i,\sigma\left(  i\right)  }\\\text{(by (\ref{sol.ps4.6.i.buv}%
), applied to}\\u=i\text{ and }v=\sigma\left(  i\right)  \text{)}%
}}+a_{k,\sigma\left(  k\right)  }^{\prime}\cdot\prod_{\substack{i\in\left\{
1,2,\ldots,n\right\}  ;\\i\neq k}}\underbrace{b_{i,\sigma\left(  i\right)  }%
}_{\substack{=a_{i,\sigma\left(  i\right)  }^{\prime}\\\text{(by
(\ref{sol.ps4.6.i.a'uv}), applied to}\\u=i\text{ and }v=\sigma\left(
i\right)  \text{)}}}\\
&  =\underbrace{a_{k,\sigma\left(  k\right)  }\cdot\prod_{\substack{i\in
\left\{  1,2,\ldots,n\right\}  ;\\i\neq k}}a_{i,\sigma\left(  i\right)  }%
}_{\substack{=\prod_{i=1}^{n}a_{i,\sigma\left(  i\right)  }\\\text{(by
(\ref{sol.ps4.6.i.sigma.pf.1}))}}}+\underbrace{a_{k,\sigma\left(  k\right)
}^{\prime}\cdot\prod_{\substack{i\in\left\{  1,2,\ldots,n\right\}  ;\\i\neq
k}}a_{i,\sigma\left(  i\right)  }^{\prime}}_{\substack{=\prod_{i=1}%
^{n}a_{i,\sigma\left(  i\right)  }^{\prime}\\\text{(by
(\ref{sol.ps4.6.i.sigma.pf.2}))}}}\\
&  =\prod_{i=1}^{n}a_{i,\sigma\left(  i\right)  }+\prod_{i=1}^{n}%
a_{i,\sigma\left(  i\right)  }^{\prime}.
\end{align*}
This proves (\ref{sol.ps4.6.i.sigma}).}.
\end{verlong}

We have $A=\left(  a_{i,j}^{\prime}\right)  _{1\leq i\leq n,\ 1\leq j\leq n}$;
therefore, (\ref{eq.det.eq.2}) (applied to $A^{\prime}$ and $a_{i,j}^{\prime}$
instead of $A$ and $a_{i,j}$) yields%
\begin{equation}
\det A^{\prime}=\sum_{\sigma\in S_{n}}\left(  -1\right)  ^{\sigma}\prod
_{i=1}^{n}a_{i,\sigma\left(  i\right)  }^{\prime}. \label{sol.ps4.6.i.detA'}%
\end{equation}

\begin{vershort}
Now, recall that $B=\left(  b_{i,j}\right)  _{1\leq i\leq n,\ 1\leq j\leq n}$.
Hence, (\ref{eq.det.eq.2}) (applied to $B$ and $b_{i,j}$ instead of $A$ and
$a_{i,j}$) yields%
\begin{align*}
\det B  &  =\sum_{\sigma\in S_{n}}\left(  -1\right)  ^{\sigma}%
\underbrace{\prod_{i=1}^{n}b_{i,\sigma\left(  i\right)  }}_{\substack{=\prod
_{i=1}^{n}a_{i,\sigma\left(  i\right)  }+\prod_{i=1}^{n}a_{i,\sigma\left(
i\right)  }^{\prime}\\\text{(by (\ref{sol.ps4.6.i.short.sigma}))}}%
}=\sum_{\sigma\in S_{n}}\left(  -1\right)  ^{\sigma}\left(  \prod_{i=1}%
^{n}a_{i,\sigma\left(  i\right)  }+\prod_{i=1}^{n}a_{i,\sigma\left(  i\right)
}^{\prime}\right) \\
&  =\underbrace{\sum_{\sigma\in S_{n}}\left(  -1\right)  ^{\sigma}\prod
_{i=1}^{n}a_{i,\sigma\left(  i\right)  }}_{\substack{=\det A\\\text{(by
(\ref{eq.det.eq.2}))}}}+\underbrace{\sum_{\sigma\in S_{n}}\left(  -1\right)
^{\sigma}\prod_{i=1}^{n}a_{i,\sigma\left(  i\right)  }^{\prime}}%
_{\substack{=\det A^{\prime}\\\text{(by (\ref{sol.ps4.6.i.detA'}))}}}=\det
A+\det A^{\prime}.
\end{align*}
This solves Exercise \ref{exe.ps4.6} \textbf{(i)}.
\end{vershort}

\begin{verlong}
Now, recall that $B=\left(  b_{i,j}\right)  _{1\leq i\leq n,\ 1\leq j\leq n}$.
Hence, (\ref{eq.det.eq.2}) (applied to $B$ and $b_{i,j}$ instead of $A$ and
$a_{i,j}$) yields%
\begin{align*}
\det B  &  =\sum_{\sigma\in S_{n}}\left(  -1\right)  ^{\sigma}%
\underbrace{\prod_{i=1}^{n}b_{i,\sigma\left(  i\right)  }}_{\substack{=\prod
_{i=1}^{n}a_{i,\sigma\left(  i\right)  }+\prod_{i=1}^{n}a_{i,\sigma\left(
i\right)  }^{\prime}\\\text{(by (\ref{sol.ps4.6.i.sigma}))}}}=\sum_{\sigma\in
S_{n}}\left(  -1\right)  ^{\sigma}\left(  \prod_{i=1}^{n}a_{i,\sigma\left(
i\right)  }+\prod_{i=1}^{n}a_{i,\sigma\left(  i\right)  }^{\prime}\right) \\
&  =\underbrace{\sum_{\sigma\in S_{n}}\left(  -1\right)  ^{\sigma}\prod
_{i=1}^{n}a_{i,\sigma\left(  i\right)  }}_{\substack{=\det A\\\text{(by
(\ref{eq.det.eq.2}))}}}+\underbrace{\sum_{\sigma\in S_{n}}\left(  -1\right)
^{\sigma}\prod_{i=1}^{n}a_{i,\sigma\left(  i\right)  }^{\prime}}%
_{\substack{=\det A^{\prime}\\\text{(by (\ref{sol.ps4.6.i.detA'}))}}}=\det
A+\det A^{\prime}.
\end{align*}
This solves Exercise \ref{exe.ps4.6} \textbf{(i)}.
\end{verlong}

\textbf{(j)} We know that $A^{\prime}$ is an $n\times n$-matrix whose columns
equal the corresponding columns of $A$ except (perhaps) the $k$-th column.
Thus, $\left(  A^{\prime}\right)  ^{T}$ is an $n\times n$-matrix whose rows
equal the corresponding rows of $A^{T}$ except (perhaps) the $k$-th row
(because the columns of $A$ correspond to the rows of $A^{T}$%
\ \ \ \ \footnote{See Remark \ref{rmk.colsA=rowsATT.rmk} for the meaning of
\textquotedblleft correspond\textquotedblright\ we are using here.}, and
similarly the columns of $A^{\prime}$ correspond to the rows of $\left(
A^{\prime}\right)  ^{T}$).

Also, we know that $B$ is the $n\times n$-matrix obtained from $A$ by adding
the $k$-th column of $A^{\prime}$ to the $k$-th column of $A$. Hence, $B^{T}$
is the $n\times n$-matrix obtained from $A^{T}$ by adding the $k$-th row of
$\left(  A^{\prime}\right)  ^{T}$ to the $k$-th row of $A^{T}$ (because the
columns of $A$ correspond to the rows of $A^{T}$, and similarly the columns of
$A^{\prime}$ correspond to the rows of $\left(  A^{\prime}\right)  ^{T}$, and
similarly the columns of $B$ correspond to the rows of $B^{T}$).

Hence, Exercise \ref{exe.ps4.6} \textbf{(i)} (applied to $A^{T}$, $\left(
A^{\prime}\right)  ^{T}$ and $B^{T}$ instead of $A$, $A^{\prime}$ and $B$)
yields $\det\left(  B^{T}\right)  =\det\left(  A^{T}\right)  +\det\left(
\left(  A^{\prime}\right)  ^{T}\right)  $.

But Exercise \ref{exe.ps4.4} yields $\det\left(  A^{T}\right)  =\det A$. Also,
Exercise \ref{exe.ps4.4} (applied to $B$ instead of $A$) yields $\det\left(
B^{T}\right)  =\det B$. Finally, Exercise \ref{exe.ps4.4} (applied to
$A^{\prime}$ instead of $A$) yields $\det\left(  \left(  A^{\prime}\right)
^{T}\right)  =\det A^{\prime}$.

But recall that $\det\left(  B^{T}\right)  =\det\left(  A^{T}\right)
+\det\left(  \left(  A^{\prime}\right)  ^{T}\right)  $. This rewrites as $\det
B=\det A+\det A^{\prime}$ (since $\det\left(  B^{T}\right)  =\det B$ and
$\det\left(  A^{T}\right)  =\det A$ and $\det\left(  \left(  A^{\prime
}\right)  ^{T}\right)  =\det A^{\prime}$). This solves Exercise
\ref{exe.ps4.6} \textbf{(j)}.

[\textit{Remark:} It is tempting to regard Exercise \ref{exe.ps4.6}
\textbf{(e)} as a consequence of Exercise \ref{exe.ps4.6} \textbf{(a)},
because if a matrix $A$ has two equal rows, then swapping these two rows does
not change the matrix $A$, and thus Exercise \ref{exe.ps4.6} \textbf{(a)}
(applied to $B=A$) yields that $\det A=-\det A$ in this case. However, we
cannot conclude $\det A=0$ from $\det A=-\det A$ in general, unless we know
that we can \textquotedblleft divide by $2$\textquotedblright\ (or at least
cancel a factor of $2$ from equalities) in the commutative ring $\mathbb{K}$.
For example, if $\mathbb{K}$ is the ring $\mathbb{Z}/2\mathbb{Z}$, then every
$a\in\mathbb{K}$ satisfies $a=-a$, but not every $a\in\mathbb{K}$ satisfies
$a=0$. Of course, if $\mathbb{K}=\mathbb{Z}$ or $\mathbb{K}=\mathbb{Q}$ or
$\mathbb{K}=\mathbb{R}$ or $\mathbb{K}=\mathbb{C}$, then this argument works
and thus Exercise \ref{exe.ps4.6} \textbf{(e)} does follow from Exercise
\ref{exe.ps4.6} \textbf{(a)} in these cases.]
\end{proof}

\subsection{Solution to Exercise \ref{exe.ps4.6k}}

Our solution to Exercise \ref{exe.ps4.6k} relies on Lemma \ref{lem.det.sigma}
and on Lemma \ref{lem.colsA=rowsATT}. Thus, we advise the reader to read the
proofs of these two lemmas before the solution of the exercise.

Before we start solving Exercise \ref{exe.ps4.6k}, let us prove a simple fact:

\begin{lemma}
\label{lem.sol.ps4.6k.lem.1}Let $n\in\mathbb{N}$. Let $A$ be an $n\times
n$-matrix such that some row of the matrix $A$ equals a scalar multiple of
some other row of $A$. Then,%
\begin{equation}
\det A=0. \label{sol.ps4.6k.lem.1}%
\end{equation}

\end{lemma}

\begin{proof}
[Proof of Lemma \ref{lem.sol.ps4.6k.lem.1}.]Some row of the matrix $A$ equals
a scalar multiple of some other row of $A$. In other words, there exist two
distinct elements $k$ and $\ell$ of $\left\{  1,2,\ldots,n\right\}  $ such
that the $k$-th row of the matrix $A$ equals a scalar multiple of the $\ell
$-th row of $A$. Consider these $k$ and $\ell$.

\begin{vershort}
Let $C$ be the $n\times n$-matrix obtained from $A$ by replacing the $k$-th
row of $A$ by the $\ell$-th row of $A$. Then, the matrix $C$ has two equal
rows (namely, its $k$-th row and its $\ell$-th row). Hence, $\det C=0$ (by
Exercise \ref{exe.ps4.6} \textbf{(e)}, applied to $C$ instead of $A$).

Recall that the $k$-th row of the matrix $A$ equals a scalar multiple of the
$\ell$-th row of $A$. In other words, there exists a $\lambda\in\mathbb{K}$
such that%
\begin{equation}
\left(  \text{the }k\text{-th row of }A\right)  =\lambda\left(  \text{the
}\ell\text{-th row of }A\right)  . \label{sol.ps4.6k.short.lem.1.pf.1}%
\end{equation}
Consider this $\lambda$. By the construction of $C$, we have
\begin{equation}
\left(  \text{the }k\text{-th row of }C\right)  =\left(  \text{the }%
\ell\text{-th row of }A\right)  . \label{sol.ps4.6k.short.lem.1.pf.2}%
\end{equation}
Thus, (\ref{sol.ps4.6k.short.lem.1.pf.1}) becomes%
\begin{align}
\left(  \text{the }k\text{-th row of }A\right)   &  =\lambda
\underbrace{\left(  \text{the }\ell\text{-th row of }A\right)  }%
_{\substack{=\left(  \text{the }k\text{-th row of }C\right)  \\\text{(by
(\ref{sol.ps4.6k.short.lem.1.pf.2}))}}}\nonumber\\
&  =\lambda\left(  \text{the }k\text{-th row of }C\right)  .
\label{sol.ps4.6k.short.lem.1.pf.3}%
\end{align}
On the other hand, for every $u\in\left\{  1,2,\ldots,n\right\}  $ satisfying
$u\neq k$, we have%
\[
\left(  \text{the }u\text{-th row of }C\right)  =\left(  \text{the }u\text{-th
row of }A\right)
\]
(since the construction of $C$ involves modifying only the $k$-th row of $A$)
and thus%
\begin{equation}
\left(  \text{the }u\text{-th row of }A\right)  =\left(  \text{the }u\text{-th
row of }C\right)  . \label{sol.ps4.6k.short.lem.1.pf.4}%
\end{equation}

Now, combining (\ref{sol.ps4.6k.short.lem.1.pf.3}) with
(\ref{sol.ps4.6k.short.lem.1.pf.4}), we obtain the following description of
$A$ from $C$: The matrix $A$ is the matrix obtained from $C$ by multiplying
the $k$-th row by $\lambda$. Then, $\det A=\lambda\det C$ (by Exercise
\ref{exe.ps4.6} \textbf{(g)}, applied to $C$ and $A$ instead of $A$ and $B$).
Thus, $\det A=\lambda\underbrace{\det C}_{=0}=0$. This proves Lemma
\ref{lem.sol.ps4.6k.lem.1}. \qedhere

\end{vershort}

\begin{verlong}
Let $w$ be the $\ell$-th row of $A$ (regarded, as usual, as a row vector).
Thus, $w=\left(  \text{the }\ell\text{-th row of }A\right)  $.

Recall that the $k$-th row of $A$ equals a scalar multiple of the $\ell$-th
row of $A$. In other words, the $k$-th row of $A$ equals a scalar multiple of
$w$ (since $w$ is the $\ell$-th row of $A$). In other words, there exists some
$\lambda\in\mathbb{K}$ such that
\begin{equation}
\left(  \text{the }k\text{-th row of }A\right)  =\lambda w.
\label{sol.ps4.6k.lem.2}%
\end{equation}
Consider this $\lambda$.

Let $C$ be the $n\times n$-matrix obtained from $A$ by replacing the $k$-th
row of $A$ by the row vector $w$. Thus,%
\begin{align}
&  \left(  \left(  \text{the }u\text{-th row of }C\right)  =\left(  \text{the
}u\text{-th row of }A\right)  \right. \label{sol.ps4.6k.lem.3}\\
&  \ \ \ \ \ \ \ \ \ \ \left.  \text{for all }u\in\left\{  1,2,\ldots
,n\right\}  \text{ satisfying }u\neq k\right)  ,\nonumber
\end{align}
whereas%
\begin{equation}
\left(  \text{the }k\text{-th row of }C\right)  =w. \label{sol.ps4.6k.lem.4}%
\end{equation}
The matrix $C$ has two equal rows\footnote{\textit{Proof.} The numbers $k$ and
$\ell$ are distinct. Thus, $\ell\neq k$. Hence, (\ref{sol.ps4.6k.lem.3})
(applied to $u=\ell$) yields%
\begin{align*}
\left(  \text{the }\ell\text{-th row of }C\right)   &  =\left(  \text{the
}\ell\text{-th row of }A\right)  =w\ \ \ \ \ \ \ \ \ \ \left(  \text{since
}w=\left(  \text{the }\ell\text{-th row of }A\right)  \right) \\
&  =\left(  \text{the }k\text{-th row of }C\right)
\ \ \ \ \ \ \ \ \ \ \left(  \text{by (\ref{sol.ps4.6k.lem.4})}\right)  .
\end{align*}
In other words, the $\ell$-th row of $C$ and the $k$-th row of $C$ are equal.
Since $\ell\neq k$, this shows that the matrix $C$ has two equal rows. Qed.}.
Hence, $\det C=0$ (by Exercise \ref{exe.ps4.6} \textbf{(e)}, applied to $C$
instead of $A$).

On the other hand, let $C^{\prime}$ be the matrix obtained from $C$ by
multiplying the $k$-th row by $\lambda$. Then, $\det C^{\prime}=\lambda\det C$
(by Exercise \ref{exe.ps4.6} \textbf{(g)}, applied to $C$ and $C^{\prime}$
instead of $A$ and $B$). Thus, $\det C^{\prime}=\lambda\underbrace{\det
C}_{=0}=0$.

Recall that $C^{\prime}$ is the matrix obtained from $C$ by multiplying the
$k$-th row by $\lambda$. In other words,%
\begin{align}
&  \left(  \left(  \text{the }u\text{-th row of }C^{\prime}\right)  =\left(
\text{the }u\text{-th row of }C\right)  \right. \label{sol.ps4.6k.lem.3'}\\
&  \ \ \ \ \ \ \ \ \ \ \left.  \text{for all }u\in\left\{  1,2,\ldots
,n\right\}  \text{ satisfying }u\neq k\right)  ,\nonumber
\end{align}
while%
\begin{equation}
\left(  \text{the }k\text{-th row of }C^{\prime}\right)  =\lambda\left(
\text{the }k\text{-th row of }C\right)  . \label{sol.ps4.6k.lem.4'}%
\end{equation}

On the other hand, for every $u\in\left\{  1,2,\ldots,n\right\}  $, we have%
\begin{equation}
\left(  \text{the }u\text{-th row of }A\right)  =\left(  \text{the }u\text{-th
row of }C^{\prime}\right)  \label{sol.ps4.6k.lem.5}%
\end{equation}
\footnote{\textit{Proof of (\ref{sol.ps4.6k.lem.5}):} Let $u\in\left\{
1,2,\ldots,n\right\}  $. We must prove (\ref{sol.ps4.6k.lem.5}).
\par
If $u\neq k$, then (\ref{sol.ps4.6k.lem.5}) follows from%
\begin{align*}
\left(  \text{the }u\text{-th row of }A\right)   &  =\left(  \text{the
}u\text{-th row of }C\right)  \ \ \ \ \ \ \ \ \ \ \left(  \text{by
(\ref{sol.ps4.6k.lem.3})}\right) \\
&  =\left(  \text{the }u\text{-th row of }C^{\prime}\right)
\ \ \ \ \ \ \ \ \ \ \left(  \text{by (\ref{sol.ps4.6k.lem.3'})}\right)  .
\end{align*}
Hence, if $u\neq k$, then (\ref{sol.ps4.6k.lem.5}) is proven. Therefore, for
the rest of this proof of (\ref{sol.ps4.6k.lem.5}), we can WLOG assume that we
don't have $u\neq k$. Assume this.
\par
We have $u=k$ (since we don't have $u\neq k$). Hence,%
\begin{align*}
\left(  \text{the }\underbrace{u}_{=k}\text{-th row of }A\right)   &  =\left(
\text{the }k\text{-th row of }A\right)  =\lambda\underbrace{w}%
_{\substack{=\left(  \text{the }k\text{-th row of }C\right)  \\\text{(by
(\ref{sol.ps4.6k.lem.4}))}}}\ \ \ \ \ \ \ \ \ \ \left(  \text{by
(\ref{sol.ps4.6k.lem.2})}\right) \\
&  =\lambda\left(  \text{the }k\text{-th row of }C\right)  =\left(  \text{the
}\underbrace{k}_{=u}\text{-th row of }C^{\prime}\right)
\ \ \ \ \ \ \ \ \ \ \left(  \text{by (\ref{sol.ps4.6k.lem.4'})}\right) \\
&  =\left(  \text{the }u\text{-th row of }C^{\prime}\right)  .
\end{align*}
Hence, (\ref{sol.ps4.6k.lem.5}) is proven.}. In other words, every row of $A$
equals the corresponding row of $C^{\prime}$. Thus, the matrices $A$ and
$C^{\prime}$ must be equal. In other words, $A=C^{\prime}$, so that $\det
A=\det C^{\prime}=0$. This proves Lemma \ref{lem.sol.ps4.6k.lem.1}.
\end{verlong}
\end{proof}

Now, let us come to the actual solution of Exercise \ref{exe.ps4.6k}.

\begin{proof}
[Solution to Exercise \ref{exe.ps4.6k}.]Let us write the $n\times n$-matrix
$A$ in the form \newline$A=\left(  a_{i,j}\right)  _{1\leq i\leq n,\ 1\leq
j\leq n}$.

\textbf{(a)} We need to prove that if we add a scalar multiple of a row of $A$
to another row of $A$, then the determinant of $A$ does not change. In other
words, we need to prove that if $B$ is an $n\times n$-matrix obtained from $A$
by adding a scalar multiple of a row of $A$ to another row of $A$, then $\det
B=\det A$.

So let $B$ be an $n\times n$-matrix obtained from $A$ by adding a scalar
multiple of a row of $A$ to another row of $A$. We then need to show that
$\det B=\det A$.

We know that $B$ is an $n\times n$-matrix obtained from $A$ by adding a scalar
multiple of a row of $A$ to another row of $A$. In other words, there exist
two distinct elements $\ell$ and $k$ of $\left\{  1,2,\ldots,n\right\}  $ such
that $B$ is the $n\times n$-matrix obtained from $A$ by adding a scalar
multiple of the $\ell$-th row of $A$ to the $k$-th row of $A$. Consider these
$\ell$ and $k$.

\begin{vershort}
We know that $B$ is the $n\times n$-matrix obtained from $A$ by adding a
scalar multiple of the $\ell$-th row of $A$ to the $k$-th row of $A$. In other
words, there exists a $\lambda\in\mathbb{K}$ such that $B$ is $n\times
n$-matrix obtained from $A$ by adding $\lambda\cdot\left(  \text{the }%
\ell\text{-th row of }A\right)  $ to the $k$-th row of $A$. Consider this
$\lambda$.

Let $A^{\prime}$ be the $n\times n$-matrix obtained from $A$ by replacing the
$k$-th row of $A$ by $\lambda\cdot\left(  \text{the }\ell\text{-th row of
}A\right)  $. Then,%
\begin{align}
\left(  \text{the }k\text{-th row of }A^{\prime}\right)   &  =\lambda
\cdot\underbrace{\left(  \text{the }\ell\text{-th row of }A\right)
}_{\substack{=\left(  \text{the }\ell\text{-th row of }A^{\prime}\right)
\\\text{(since the }\ell\text{-th row of }A^{\prime}\text{ has been}%
\\\text{taken over from }A\text{ unchanged)}}}\label{sol.ps4.6k.a.short.1}\\
&  =\lambda\cdot\left(  \text{the }\ell\text{-th row of }A^{\prime}\right)
.\nonumber
\end{align}
Thus, the $k$-th row of $A^{\prime}$ equals a scalar multiple of the $\ell$-th
row of $A^{\prime}$. Thus, some row of the matrix $A^{\prime}$ equals a scalar
multiple of some other row of $A^{\prime}$. Hence, (\ref{sol.ps4.6k.lem.1})
(applied to $A^{\prime}$ instead of $A$) shows that $\det A^{\prime}=0$.

On the other hand, $A^{\prime}$ is an $n\times n$-matrix whose rows equal the
corresponding rows of $A$ except (perhaps) the $k$-th row (because of the
construction of $A^{\prime}$). Also, $B$ is the $n\times n$-matrix obtained
from $A$ by adding $\lambda\cdot\left(  \text{the }\ell\text{-th row of
}A\right)  $ to the $k$-th row of $A$. Because of (\ref{sol.ps4.6k.a.short.1}%
), this can be rewritten as follows: $B$ is the $n\times n$-matrix obtained
from $A$ by adding the $k$-th row of $A^{\prime}$ to the $k$-th row of $A$.
Exercise \ref{exe.ps4.6} \textbf{(i)} thus yields $\det B=\det
A+\underbrace{\det A^{\prime}}_{=0}=\det A$.
\end{vershort}

\begin{verlong}
We have $\ell\neq k$ (since $\ell$ and $k$ are distinct).

We know that $B$ is the $n\times n$-matrix obtained from $A$ by adding a
scalar multiple of the $\ell$-th row of $A$ to the $k$-th row of $A$. In other
words, we have%
\begin{align}
&  \left(  \left(  \text{the }u\text{-th row of }B\right)  =\left(  \text{the
}u\text{-th row of }A\right)  \right. \label{sol.ps4.6k.a.1}\\
&  \ \ \ \ \ \ \ \ \ \ \left.  \text{for all }u\in\left\{  1,2,\ldots
,n\right\}  \text{ satisfying }u\neq k\right)  ,\nonumber
\end{align}
while the $k$-th row of $B$ is the sum of the $k$-th row of $A$ with a scalar
multiple of the $\ell$-th row of $A$.

We know that the $k$-th row of $B$ is the sum of the $k$-th row of $A$ with a
scalar multiple of the $\ell$-th row of $A$. In other words,%
\begin{equation}
\left(  \text{the }k\text{-th row of }B\right)  =\left(  \text{the }k\text{-th
row of }A\right)  +w, \label{sol.ps4.6k.a.2}%
\end{equation}
where $w$ is some scalar multiple of the $\ell$-th row of $A$. Consider this
$w$. There exists a $\lambda\in\mathbb{K}$ such that
\begin{equation}
w=\lambda\left(  \text{the }\ell\text{-th row of }A\right)
\label{sol.ps4.6k.a.3}%
\end{equation}
(since $w$ is some scalar multiple of the $\ell$-th row of $A$). Consider this
$\lambda$.

Let $A^{\prime}$ be the $n\times n$-matrix obtained from $A$ by replacing the
$k$-th row of $A$ by the row vector $w$. Thus,%
\begin{align}
&  \left(  \left(  \text{the }u\text{-th row of }A^{\prime}\right)  =\left(
\text{the }u\text{-th row of }A\right)  \right. \label{sol.ps4.6k.a.5}\\
&  \ \ \ \ \ \ \ \ \ \ \left.  \text{for all }u\in\left\{  1,2,\ldots
,n\right\}  \text{ satisfying }u\neq k\right)  ,\nonumber
\end{align}
whereas%
\begin{equation}
\left(  \text{the }k\text{-th row of }A^{\prime}\right)  =w.
\label{sol.ps4.6k.a.6}%
\end{equation}
Then, the $k$-th row of the matrix $A^{\prime}$ equals a scalar multiple of
the $\ell$-th row of $A^{\prime}$\ \ \ \ \footnote{\textit{Proof.} We have
$\ell\neq k$. Thus, (\ref{sol.ps4.6k.a.5}) (applied to $u=\ell$) yields
$\left(  \text{the }\ell\text{-th row of }A^{\prime}\right)  =\left(
\text{the }\ell\text{-th row of }A\right)  $. Hence, $\lambda
\underbrace{\left(  \text{the }\ell\text{-th row of }A^{\prime}\right)
}_{=\left(  \text{the }\ell\text{-th row of }A\right)  }=\lambda\left(
\text{the }\ell\text{-th row of }A\right)  =w$ (by (\ref{sol.ps4.6k.a.3})).
Compared with (\ref{sol.ps4.6k.a.6}), this yields $\left(  \text{the
}k\text{-th row of }A^{\prime}\right)  =\lambda\left(  \text{the }%
\ell\text{-th row of }A^{\prime}\right)  $. Hence, the $k$-th row of the
matrix $A^{\prime}$ equals a scalar multiple of the $\ell$-th row of
$A^{\prime}$. Qed.}. Therefore, some row of the matrix $A^{\prime}$ equals a
scalar multiple of some other row of $A^{\prime}$. Hence,
(\ref{sol.ps4.6k.lem.1}) (applied to $A^{\prime}$ instead of $A$) yields $\det
A^{\prime}=0$.

On the other hand, $A^{\prime}$ is an $n\times n$-matrix whose rows equal the
corresponding rows of $A$ except (perhaps) the $k$-th row\footnote{In fact,
(\ref{sol.ps4.6k.a.5}) says precisely that the rows of $A^{\prime}$ equal the
corresponding rows of $A$ except (perhaps) the $k$-th row.}.

Now, let $B^{\prime}$ be the $n\times n$-matrix obtained from $A$ by adding
the $k$-th row of $A^{\prime}$ to the $k$-th row of $A$. Thus,%
\begin{align}
&  \left(  \left(  \text{the }u\text{-th row of }B^{\prime}\right)  =\left(
\text{the }u\text{-th row of }A\right)  \right. \label{sol.ps4.6k.a.11}\\
&  \ \ \ \ \ \ \ \ \ \ \left.  \text{for all }u\in\left\{  1,2,\ldots
,n\right\}  \text{ satisfying }u\neq k\right)  ,\nonumber
\end{align}
whereas%
\begin{equation}
\left(  \text{the }k\text{-th row of }B^{\prime}\right)  =\left(  \text{the
}k\text{-th row of }A\right)  +\left(  \text{the }k\text{-th row of }%
A^{\prime}\right)  . \label{sol.ps4.6k.a.12}%
\end{equation}

Exercise \ref{exe.ps4.6} \textbf{(i)} (applied to $B^{\prime}$ instead of $B$)
yields $\det B^{\prime}=\det A+\underbrace{\det A^{\prime}}_{=0}=\det A$.

Now, for every $u\in\left\{  1,2,\ldots,n\right\}  $, we have%
\begin{equation}
\left(  \text{the }u\text{-th row of }B^{\prime}\right)  =\left(  \text{the
}u\text{-th row of }B\right)  \label{sol.ps4.6k.a.13}%
\end{equation}
\footnote{\textit{Proof of (\ref{sol.ps4.6k.a.13}):} Let $u\in\left\{
1,2,\ldots,n\right\}  $. We must prove (\ref{sol.ps4.6k.a.13}).
\par
If $u\neq k$, then (\ref{sol.ps4.6k.a.13}) follows from%
\begin{align*}
\left(  \text{the }u\text{-th row of }B^{\prime}\right)   &  =\left(
\text{the }u\text{-th row of }A\right)  \ \ \ \ \ \ \ \ \ \ \left(  \text{by
(\ref{sol.ps4.6k.a.11})}\right) \\
&  =\left(  \text{the }u\text{-th row of }B\right)
\ \ \ \ \ \ \ \ \ \ \left(  \text{by (\ref{sol.ps4.6k.a.1})}\right)  .
\end{align*}
Hence, if $u\neq k$, then (\ref{sol.ps4.6k.a.13}) is proven. Therefore, for
the rest of this proof of (\ref{sol.ps4.6k.a.13}), we can WLOG assume that we
don't have $u\neq k$. Assume this.
\par
We have $u=k$ (since we don't have $u\neq k$). Hence,%
\begin{align*}
\left(  \text{the }\underbrace{u}_{=k}\text{-th row of }B^{\prime}\right)   &
=\left(  \text{the }k\text{-th row of }B^{\prime}\right) \\
&  =\left(  \text{the }k\text{-th row of }A\right)  +\underbrace{\left(
\text{the }k\text{-th row of }A^{\prime}\right)  }_{\substack{=w\\\text{(by
(\ref{sol.ps4.6k.a.6}))}}}\ \ \ \ \ \ \ \ \ \ \left(  \text{by
(\ref{sol.ps4.6k.lem.2})}\right) \\
&  =\left(  \text{the }k\text{-th row of }A\right)  +w=\left(  \text{the
}\underbrace{k}_{=u}\text{-th row of }B\right)  \ \ \ \ \ \ \ \ \ \ \left(
\text{by (\ref{sol.ps4.6k.a.2})}\right) \\
&  =\left(  \text{the }u\text{-th row of }B\right)  .
\end{align*}
Hence, (\ref{sol.ps4.6k.a.13}) is proven.}. In other words, every row of
$B^{\prime}$ equals the corresponding row of $B$. Thus, the matrices
$B^{\prime}$ and $B$ must be equal. In other words, $B^{\prime}=B$, so that
$\det B^{\prime}=\det B$.

Compared with $\det B^{\prime}=\det A$, this yields $\det B=\det A$.
\end{verlong}

Thus, we have shown that%
\begin{equation}
\det B=\det A. \label{sol.ps4.6k.a.result}%
\end{equation}
This completes our solution to Exercise \ref{exe.ps4.6k} \textbf{(a)}.

\textbf{(b)} We need to prove that if we add a scalar multiple of a column of
$A$ to another column of $A$, then the determinant of $A$ does not change. In
other words, we need to prove that if $B$ is an $n\times n$-matrix obtained
from $A$ by adding a scalar multiple of a column of $A$ to another column of
$A$, then $\det B=\det A$.

So let $B$ be an $n\times n$-matrix obtained from $A$ by adding a scalar
multiple of a column of $A$ to another column of $A$. We then need to show
that $\det B=\det A$.

We know that $B$ is an $n\times n$-matrix obtained from $A$ by adding a scalar
multiple of a column of $A$ to another column of $A$. Therefore, $B^{T}$ is an
$n\times n$-matrix obtained from $A^{T}$ by adding a scalar multiple of a row
of $A^{T}$ to another row of $A^{T}$ (because the columns of $A$ correspond to
the rows of $A^{T}$\ \ \ \ \footnote{See Remark \ref{rmk.colsA=rowsATT.rmk}
for the meaning of \textquotedblleft correspond\textquotedblright\ we are
using here.}). Therefore, (\ref{sol.ps4.6k.a.result}) (applied to $A^{T}$ and
$B^{T}$ instead of $A$ and $B$) yields $\det\left(  B^{T}\right)  =\det\left(
A^{T}\right)  $.

But Exercise \ref{exe.ps4.4} yields $\det\left(  A^{T}\right)  =\det A$. Also,
Exercise \ref{exe.ps4.4} (applied to $B$ instead of $A$) yields $\det\left(
B^{T}\right)  =\det B$. But recall that $\det\left(  B^{T}\right)
=\det\left(  A^{T}\right)  $. This rewrites as $\det B=\det A$ (since
$\det\left(  B^{T}\right)  =\det B$ and $\det\left(  A^{T}\right)  =\det A$).
This solves Exercise \ref{exe.ps4.6k} \textbf{(b)}.
\end{proof}

\subsection{Solution to Exercise \ref{exe.prodrule}}

Our goal is now to prove Lemma \ref{lem.prodrule}. Actually, we will prove a
more general result:

\begin{lemma}
\label{lem.prodrule.S}Let $n\in\mathbb{N}$. For every $i\in\left\{
1,2,\ldots,n\right\}  $, let $Z_{i}$ be a finite set. For every $i\in\left\{
1,2,\ldots,n\right\}  $ and every $k\in Z_{i}$, let $p_{i,k}$ be an element of
$\mathbb{K}$. Then,%
\[
\prod_{i=1}^{n}\ \ \sum_{k\in Z_{i}}p_{i,k}=\sum_{\left(  k_{1},k_{2}%
,\ldots,k_{n}\right)  \in Z_{1}\times Z_{2}\times\cdots\times Z_{n}}%
\ \ \prod_{i=1}^{n}p_{i,k_{i}}.
\]

\end{lemma}

Let us first show the particular case of this lemma for $n=2$:

\begin{lemma}
\label{lem.prodrule.S.n=2}Let $X$ and $Y$ be two finite sets. For every $x\in
X$, let $q_{x}$ be an element of $\mathbb{K}$. For every $y\in Y$, let $r_{y}$
be an element of $\mathbb{K}$. Then,%
\[
\left(  \sum_{x\in X}q_{x}\right)  \left(  \sum_{y\in Y}r_{y}\right)
=\sum_{\left(  x,y\right)  \in X\times Y}q_{x}r_{y}.
\]

\end{lemma}

\begin{vershort}
\begin{proof}
[Proof of Lemma \ref{lem.prodrule.S.n=2}.]We have%
\[
\underbrace{\sum_{\left(  x,y\right)  \in X\times Y}}_{=\sum_{x\in X}%
\ \ \sum_{y\in Y}}q_{x}r_{y}=\sum_{x\in X}\ \ \underbrace{\sum_{y\in Y}%
q_{x}r_{y}}_{=q_{x}\sum_{y\in Y}r_{y}}=\sum_{x\in X}\left(  q_{x}\sum_{y\in
Y}r_{y}\right)  =\left(  \sum_{x\in X}q_{x}\right)  \left(  \sum_{y\in Y}%
r_{y}\right)  .
\]
This proves Lemma \ref{lem.prodrule.S.n=2}.
\end{proof}
\end{vershort}

\begin{verlong}
\begin{proof}
[Proof of Lemma \ref{lem.prodrule.S.n=2}.]The equalities (\ref{eq.sum.linear2}%
) and (\ref{eq.sum.fubini}) (and most other properties of finite sums) remain
valid if $\mathbb{A}$ is replaced by $\mathbb{K}$ throughout them. It is in
this variant that they will be used in the following proof.

For every $\lambda\in\mathbb{K}$, we have%
\begin{align}
\sum_{y\in Y}\lambda r_{y}  &  =\sum_{s\in Y}\lambda r_{s}%
\ \ \ \ \ \ \ \ \ \ \left(  \text{here, we renamed the summation index
}y\text{ as }s\right) \nonumber\\
&  =\lambda\sum_{s\in Y}r_{s}\ \ \ \ \ \ \ \ \ \ \left(
\begin{array}
[c]{c}%
\text{by (\ref{eq.sum.linear2}) (applied to }Y\text{, }r_{s}\text{ and
}\mathbb{K}\\
\text{instead of }S\text{, }a_{s}\text{ and }\mathbb{A}\text{)}%
\end{array}
\right) \nonumber\\
&  =\lambda\sum_{y\in Y}r_{y}\ \ \ \ \ \ \ \ \ \ \left(  \text{here, we
renamed the summation index }s\text{ as }y\right) \nonumber\\
&  =\left(  \sum_{y\in Y}r_{y}\right)  \lambda.
\label{pf.lem.prodrule.S.n=2.1}%
\end{align}
For every $\lambda\in\mathbb{K}$, we have%
\begin{align}
\sum_{x\in X}\lambda q_{x}  &  =\sum_{s\in X}\lambda q_{s}%
\ \ \ \ \ \ \ \ \ \ \left(  \text{here, we renamed the summation index
}x\text{ as }s\right) \nonumber\\
&  =\lambda\sum_{s\in X}q_{s}\ \ \ \ \ \ \ \ \ \ \left(
\begin{array}
[c]{c}%
\text{by (\ref{eq.sum.linear2}) (applied to }X\text{, }q_{s}\text{ and
}\mathbb{K}\\
\text{instead of }S\text{, }a_{s}\text{ and }\mathbb{A}\text{)}%
\end{array}
\right) \nonumber\\
&  =\lambda\sum_{x\in X}q_{x}\ \ \ \ \ \ \ \ \ \ \left(  \text{here, we
renamed the summation index }s\text{ as }x\right) \nonumber\\
&  =\left(  \sum_{x\in X}q_{x}\right)  \lambda.
\label{pf.lem.prodrule.S.n=2.2}%
\end{align}

From (\ref{eq.sum.fubini}) (applied to $\mathbb{K}$ and $q_{x}r_{y}$ instead
of $\mathbb{A}$ and $a_{\left(  x,y\right)  }$), we obtain%
\[
\sum_{x\in X}\ \ \sum_{y\in Y}q_{x}r_{y}=\sum_{\left(  x,y\right)  \in X\times
Y}q_{x}r_{y}=\sum_{y\in Y}\ \ \sum_{x\in X}q_{x}r_{y}.
\]
Hence,%
\begin{align*}
\sum_{\left(  x,y\right)  \in X\times Y}q_{x}r_{y}  &  =\sum_{x\in
X}\ \ \underbrace{\sum_{y\in Y}q_{x}r_{y}}_{\substack{=\left(  \sum_{y\in
Y}r_{y}\right)  q_{x}\\\text{(by (\ref{pf.lem.prodrule.S.n=2.1}) (applied to
}\lambda=q_{x}\text{))}}}=\sum_{x\in X}\left(  \sum_{y\in Y}r_{y}\right)
q_{x}\\
&  =\left(  \sum_{x\in X}q_{x}\right)  \left(  \sum_{y\in Y}r_{y}\right)
\end{align*}
(by (\ref{pf.lem.prodrule.S.n=2.2}) (applied to $\lambda=\sum_{y\in Y}r_{y}%
$)). This proves Lemma \ref{lem.prodrule.S.n=2}.
\end{proof}
\end{verlong}

\begin{proof}
[Proof of Lemma \ref{lem.prodrule.S}.]We claim that%
\begin{equation}
\prod_{i=1}^{m}\ \ \sum_{k\in Z_{i}}p_{i,k}=\sum_{\left(  k_{1},k_{2}%
,\ldots,k_{m}\right)  \in Z_{1}\times Z_{2}\times\cdots\times Z_{m}}%
\ \ \prod_{i=1}^{m}p_{i,k_{i}} \label{pf.lem.prodrule.S.claim}%
\end{equation}
for every $m\in\left\{  0,1,\ldots,n\right\}  $.

[\textit{Proof of (\ref{pf.lem.prodrule.S.claim}):} We shall prove
(\ref{pf.lem.prodrule.S.claim}) by induction over $m$:

\textit{Induction base:} Comparing
\[
\prod_{i=1}^{0}\ \ \sum_{k\in Z_{i}}p_{i,k}=\left(  \text{empty product}%
\right)  =1
\]
with%
\begin{align*}
\sum_{\left(  k_{1},k_{2},\ldots,k_{0}\right)  \in Z_{1}\times Z_{2}%
\times\cdots\times Z_{0}}\ \ \underbrace{\prod_{i=1}^{0}p_{i,k_{i}}}_{=\left(
\text{empty product}\right)  =1}  &  =\sum_{\left(  k_{1},k_{2},\ldots
,k_{0}\right)  \in Z_{1}\times Z_{2}\times\cdots\times Z_{0}}1\\
&  =\underbrace{\left\vert Z_{1}\times Z_{2}\times\cdots\times Z_{0}%
\right\vert }_{\substack{=1\\\text{(since }Z_{1}\times Z_{2}\times\cdots\times
Z_{0}\text{ is an}\\\text{empty Cartesian product)}}}\cdot1=1,
\end{align*}
we obtain $\prod_{i=1}^{0}\ \ \sum_{k\in Z_{i}}p_{i,k}=\sum_{\left(
k_{1},k_{2},\ldots,k_{0}\right)  \in Z_{1}\times Z_{2}\times\cdots\times
Z_{0}}\ \ \prod_{i=1}^{0}p_{i,k_{i}}$. In other words,
(\ref{pf.lem.prodrule.S.claim}) holds for $m=0$. This completes the induction base.

\textit{Induction step:} Let $M\in\left\{  0,1,\ldots,n\right\}  $ be
positive. Assume that (\ref{pf.lem.prodrule.S.claim}) holds for $m=M-1$. We
now must show that (\ref{pf.lem.prodrule.S.claim}) holds for $m=M$.

We have assumed that (\ref{pf.lem.prodrule.S.claim}) holds for $m=M-1$. In
other words, we have%
\begin{equation}
\prod_{i=1}^{M-1}\ \ \sum_{k\in Z_{i}}p_{i,k}=\sum_{\left(  k_{1},k_{2}%
,\ldots,k_{M-1}\right)  \in Z_{1}\times Z_{2}\times\cdots\times Z_{M-1}%
}\ \ \prod_{i=1}^{M-1}p_{i,k_{i}}. \label{pf.lem.prodrule.S.indhyp}%
\end{equation}

Lemma \ref{lem.prodrule.prod-assM} shows that the map
\begin{align*}
Z_{1}\times Z_{2}\times\cdots\times Z_{M}  &  \rightarrow\left(  Z_{1}\times
Z_{2}\times\cdots\times Z_{M-1}\right)  \times Z_{M},\\
\left(  s_{1},s_{2},\ldots,s_{M}\right)   &  \mapsto\left(  \left(
s_{1},s_{2},\ldots,s_{M-1}\right)  ,s_{M}\right)
\end{align*}
is a bijection.

For every $\left(  k_{1},k_{2},\ldots,k_{M-1}\right)  \in Z_{1}\times
Z_{2}\times\cdots\times Z_{M-1}$, we define an element \newline$g_{\left(
k_{1},k_{2},\ldots,k_{M-1}\right)  }\in\mathbb{K}$ by
\begin{equation}
g_{\left(  k_{1},k_{2},\ldots,k_{M-1}\right)  }=\prod_{i=1}^{M-1}p_{i,k_{i}}.
\label{pf.lem.prodrule.S.g=}%
\end{equation}

Now, (\ref{pf.lem.prodrule.S.indhyp}) becomes%
\begin{align}
\prod_{i=1}^{M-1}\ \ \sum_{k\in Z_{i}}p_{i,k}  &  =\sum_{\left(  k_{1}%
,k_{2},\ldots,k_{M-1}\right)  \in Z_{1}\times Z_{2}\times\cdots\times Z_{M-1}%
}\ \ \underbrace{\prod_{i=1}^{M-1}p_{i,k_{i}}}_{\substack{=g_{\left(
k_{1},k_{2},\ldots,k_{M-1}\right)  }\\\text{(by (\ref{pf.lem.prodrule.S.g=}%
))}}}\nonumber\\
&  =\sum_{\left(  k_{1},k_{2},\ldots,k_{M-1}\right)  \in Z_{1}\times
Z_{2}\times\cdots\times Z_{M-1}}g_{\left(  k_{1},k_{2},\ldots,k_{M-1}\right)
}\nonumber\\
&  =\sum_{x\in Z_{1}\times Z_{2}\times\cdots\times Z_{M-1}}g_{x}
\label{pf.lem.prodrule.S.L1}%
\end{align}
(here, we have renamed the summation index $\left(  k_{1},k_{2},\ldots
,k_{M-1}\right)  $ as $x$).

The sets $Z_{1},Z_{2},\ldots,Z_{M-1}$ are finite. Hence, their Cartesian
product $Z_{1}\times Z_{2}\times\cdots\times Z_{M-1}$ is also finite. Every
$\left(  s_{1},s_{2},\ldots,s_{M}\right)  \in Z_{1}\times Z_{2}\times
\cdots\times Z_{M}$ satisfies%
\begin{equation}
\prod_{i=1}^{M}p_{i,s_{i}}=\left(  \prod_{i=1}^{M-1}p_{i,s_{i}}\right)
p_{M,s_{M}} \label{pf.lem.prodrule.S.splioff}%
\end{equation}
(indeed, this follows by splitting off the factor for $i=M$ from the product
$\prod_{i=1}^{M}p_{i,s_{i}}$).

Now, if we split off the factor for $i=M$ from the product $\prod_{i=1}%
^{M}\ \ \sum_{k\in Z_{i}}p_{i,k}$, we obtain
\begin{align*}
&  \prod_{i=1}^{M}\ \ \sum_{k\in Z_{i}}p_{i,k}\\
&  =\underbrace{\left(  \prod_{i=1}^{M-1}\ \ \sum_{k\in Z_{i}}p_{i,k}\right)
}_{\substack{=\sum_{x\in Z_{1}\times Z_{2}\times\cdots\times Z_{M-1}}%
g_{x}\\\text{(by (\ref{pf.lem.prodrule.S.L1}))}}}\underbrace{\left(
\sum_{k\in Z_{M}}p_{M,k}\right)  }_{\substack{=\sum_{y\in Z_{M}}%
p_{M,y}\\\text{(here, we renamed the}\\\text{summation index }k\text{ as
}y\text{)}}}\\
&  =\left(  \sum_{x\in Z_{1}\times Z_{2}\times\cdots\times Z_{M-1}}%
g_{x}\right)  \left(  \sum_{y\in Z_{M}}p_{M,y}\right)  =\sum_{\left(
x,y\right)  \in\left(  Z_{1}\times Z_{2}\times\cdots\times Z_{M-1}\right)
\times Z_{M}}g_{x}p_{M,y}\\
&  \ \ \ \ \ \ \ \ \ \ \ \ \ \ \ \ \ \ \ \ \left(
\begin{array}
[c]{c}%
\text{by Lemma \ref{lem.prodrule.S.n=2}, applied to}\\
X=Z_{1}\times Z_{2}\times\cdots\times Z_{M-1}\text{, }Y=Z_{M}\text{, }%
q_{x}=g_{x}\text{ and }r_{y}=p_{M,y}%
\end{array}
\right) \\
&  =\sum_{\left(  s_{1},s_{2},\ldots,s_{M}\right)  \in Z_{1}\times Z_{2}%
\times\cdots\times Z_{M}}\underbrace{g_{\left(  s_{1},s_{2},\ldots
,s_{M-1}\right)  }}_{\substack{=\prod_{i=1}^{M-1}p_{i,s_{i}}\\\text{(by the
definition of }g_{\left(  s_{1},s_{2},\ldots,s_{M-1}\right)  }\text{)}%
}}p_{M,s_{M}}\\
&  \ \ \ \ \ \ \ \ \ \ \ \ \ \ \ \ \ \ \ \ \left(
\begin{array}
[c]{c}%
\text{here, we have substituted }\left(  \left(  s_{1},s_{2},\ldots
,s_{M-1}\right)  ,s_{M}\right)  \text{ for }\left(  x,y\right) \\
\text{in the sum, since the map}\\
Z_{1}\times Z_{2}\times\cdots\times Z_{M}\rightarrow\left(  Z_{1}\times
Z_{2}\times\cdots\times Z_{M-1}\right)  \times Z_{M},\\
\left(  s_{1},s_{2},\ldots,s_{M}\right)  \mapsto\left(  \left(  s_{1}%
,s_{2},\ldots,s_{M-1}\right)  ,s_{M}\right) \\
\text{is a bijection}%
\end{array}
\right) \\
&  =\sum_{\left(  s_{1},s_{2},\ldots,s_{M}\right)  \in Z_{1}\times Z_{2}%
\times\cdots\times Z_{M}}\underbrace{\left(  \prod_{i=1}^{M-1}p_{i,s_{i}%
}\right)  p_{M,s_{M}}}_{\substack{=\prod_{i=1}^{M}p_{i,s_{i}}\\\text{(by
(\ref{pf.lem.prodrule.S.splioff}))}}}\\
&  =\sum_{\left(  s_{1},s_{2},\ldots,s_{M}\right)  \in Z_{1}\times Z_{2}%
\times\cdots\times Z_{M}}\ \ \prod_{i=1}^{M}p_{i,s_{i}}=\sum_{\left(
k_{1},k_{2},\ldots,k_{M}\right)  \in Z_{1}\times Z_{2}\times\cdots\times
Z_{M}}\ \ \prod_{i=1}^{M}p_{i,k_{i}}%
\end{align*}
(here, we have renamed the summation index $\left(  s_{1},s_{2},\ldots
,s_{M}\right)  $ as $\left(  k_{1},k_{2},\ldots,k_{M}\right)  $). In other
words, (\ref{pf.lem.prodrule.S.claim}) holds for $m=M$. This completes the
induction step. Thus, (\ref{pf.lem.prodrule.S.claim}) is proven by induction.]

Now, (\ref{pf.lem.prodrule.S.claim}) (applied to $m=n$) yields%
\[
\prod_{i=1}^{n}\ \ \sum_{k\in Z_{i}}p_{i,k}=\sum_{\left(  k_{1},k_{2}%
,\ldots,k_{n}\right)  \in Z_{1}\times Z_{2}\times\cdots\times Z_{n}}%
\ \ \prod_{i=1}^{n}p_{i,k_{i}}.
\]
This proves Lemma \ref{lem.prodrule.S}.
\end{proof}

\begin{proof}
[Proof of Lemma \ref{lem.prodrule}.]For every $i\in\left[  n\right]  $, the
set $\left[  m_{i}\right]  $ is clearly a finite set. For every $i\in\left[
n\right]  $, we know that $p_{i,1},p_{i,2},\ldots,p_{i,m_{i}}$ are elements of
$\mathbb{K}$. In other words, for every $i\in\left[  n\right]  $ and every
$k\in\left[  m_{i}\right]  $, we know that $p_{i,k}$ is an element of
$\mathbb{K}$. In other words, for every $i\in\left\{  1,2,\ldots,n\right\}  $
and every $k\in\left[  m_{i}\right]  $, we know that $p_{i,k}$ is an element
of $\mathbb{K}$ (since $\left[  n\right]  =\left\{  1,2,\ldots,n\right\}  $).
Hence, Lemma \ref{lem.prodrule.S} (applied to $Z_{i}=\left[  m_{i}\right]  $)
yields%
\[
\prod_{i=1}^{n}\ \ \sum_{k\in\left[  m_{i}\right]  }p_{i,k}=\sum_{\left(
k_{1},k_{2},\ldots,k_{n}\right)  \in\left[  m_{1}\right]  \times\left[
m_{2}\right]  \times\cdots\times\left[  m_{n}\right]  }\ \ \prod_{i=1}%
^{n}p_{i,k_{i}}.
\]
Thus,%
\[
\sum_{\left(  k_{1},k_{2},\ldots,k_{n}\right)  \in\left[  m_{1}\right]
\times\left[  m_{2}\right]  \times\cdots\times\left[  m_{n}\right]  }%
\ \ \prod_{i=1}^{n}p_{i,k_{i}}=\prod_{i=1}^{n}\ \ \underbrace{\sum
_{k\in\left[  m_{i}\right]  }}_{=\sum_{k=1}^{m_{i}}}p_{i,k}=\prod_{i=1}%
^{n}\ \ \sum_{k=1}^{m_{i}}p_{i,k}.
\]
This proves Lemma \ref{lem.prodrule}.
\end{proof}

We have now proven Lemma \ref{lem.prodrule}, and thus solved Exercise
\ref{exe.prodrule}.

We note that Lemma \ref{lem.prodrule.S} generalizes the well-known fact that
any $n\in\mathbb{N}$ and any $n$ finite sets $A_{1},A_{2},\ldots,A_{n}$
satisfy%
\[
\left\vert A_{1}\times A_{2}\times\cdots\times A_{n}\right\vert =\left\vert
A_{1}\right\vert \cdot\left\vert A_{2}\right\vert \cdot\cdots\cdot\left\vert
A_{n}\right\vert .
\]
Indeed, if $n\in\mathbb{N}$ and if $A_{1},A_{2},\ldots,A_{n}$ are $n$ finite
sets, then Lemma \ref{lem.prodrule.S} (applied to $Z_{i}=A_{i}$ and
$p_{i,k}=1$) yields%
\begin{align*}
\prod_{i=1}^{n}\ \ \sum_{k\in A_{i}}1  &  =\sum_{\left(  k_{1},k_{2}%
,\ldots,k_{n}\right)  \in A_{1}\times A_{2}\times\cdots\times A_{n}%
}\ \ \underbrace{\prod_{i=1}^{n}1}_{=1}=\sum_{\left(  k_{1},k_{2},\ldots
,k_{n}\right)  \in A_{1}\times A_{2}\times\cdots\times A_{n}}1\\
&  =\left\vert A_{1}\times A_{2}\times\cdots\times A_{n}\right\vert
\cdot1=\left\vert A_{1}\times A_{2}\times\cdots\times A_{n}\right\vert
\end{align*}
and thus%
\[
\left\vert A_{1}\times A_{2}\times\cdots\times A_{n}\right\vert =\prod
_{i=1}^{n}\ \ \underbrace{\sum_{k\in A_{i}}1}_{=\left\vert A_{i}\right\vert
\cdot1=\left\vert A_{i}\right\vert }=\prod_{i=1}^{n}\left\vert A_{i}%
\right\vert =\left\vert A_{1}\right\vert \cdot\left\vert A_{2}\right\vert
\cdot\cdots\cdot\left\vert A_{n}\right\vert .
\]

\subsection{Solution to Exercise \ref{exe.ps4.det.fibo}}

We shall solve Exercise \ref{exe.ps4.det.fibo} more or less as you would
expect, by modifying our proof of (\ref{eq.exam.det(AB).fibo.1}) in some places.

\begin{proof}
[Solution to Exercise \ref{exe.ps4.det.fibo}.]Let $B$ be the $2\times2$-matrix
$\left(
\begin{array}
[c]{cc}%
a & 1\\
b & 0
\end{array}
\right)  $. Thus,%
\[
\det B=\det\left(
\begin{array}
[c]{cc}%
a & 1\\
b & 0
\end{array}
\right)  =a\cdot0-1\cdot b=-b.
\]

Let $A$ be the $2\times2$-matrix $\left(
\begin{array}
[c]{cc}%
x_{k+2} & x_{k+1}\\
x_{1} & x_{0}%
\end{array}
\right)  $. Then, $\det A=\det\left(
\begin{array}
[c]{cc}%
x_{k+2} & x_{k+1}\\
x_{1} & x_{0}%
\end{array}
\right)  =x_{k+2}x_{0}-x_{k+1}x_{1}$.

We now claim that%
\begin{equation}
AB^{m}=\left(
\begin{array}
[c]{cc}%
x_{m+k+2} & x_{m+k+1}\\
x_{m+1} & x_{m}%
\end{array}
\right)  \ \ \ \ \ \ \ \ \ \ \text{for every }m\in\mathbb{N}.
\label{sol.ps4.det.fibo.BmC}%
\end{equation}

[\textit{Proof of (\ref{sol.ps4.det.fibo.BmC}):} We shall prove
(\ref{sol.ps4.det.fibo.BmC}) by induction over $m$:

\textit{Induction base:} We have $A\underbrace{B^{0}}_{=I_{2}}=AI_{2}%
=A=\left(
\begin{array}
[c]{cc}%
x_{k+2} & x_{k+1}\\
x_{1} & x_{0}%
\end{array}
\right)  $. Compared with $\left(
\begin{array}
[c]{cc}%
x_{0+k+2} & x_{0+k+1}\\
x_{0+1} & x_{0}%
\end{array}
\right)  =\left(
\begin{array}
[c]{cc}%
x_{k+2} & x_{k+1}\\
x_{1} & x_{0}%
\end{array}
\right)  $, this yields $AB^{0}=\left(
\begin{array}
[c]{cc}%
x_{0+k+2} & x_{0+k+1}\\
x_{0+1} & x_{0}%
\end{array}
\right)  $. In other words, (\ref{sol.ps4.det.fibo.BmC}) holds for $m=0$. This
completes the induction base.

\textit{Induction step:} Let $M$ be a positive integer. Assume that
(\ref{sol.ps4.det.fibo.BmC}) holds for $m=M-1$. We need to show that
(\ref{sol.ps4.det.fibo.BmC}) holds for $m=M$.

We have assumed that (\ref{sol.ps4.det.fibo.BmC}) holds for $m=M-1$. In other
words,%
\[
AB^{M-1}=\left(
\begin{array}
[c]{cc}%
x_{\left(  M-1\right)  +k+2} & x_{\left(  M-1\right)  +k+1}\\
x_{\left(  M-1\right)  +1} & x_{M-1}%
\end{array}
\right)  =\left(
\begin{array}
[c]{cc}%
x_{M+k+1} & x_{M+k}\\
x_{M} & x_{M-1}%
\end{array}
\right)  .
\]
Now,%
\begin{align}
A\underbrace{B^{M}}_{=B^{M-1}\cdot B}  &  =\underbrace{AB^{M-1}}_{=\left(
\begin{array}
[c]{cc}%
x_{M+k+1} & x_{M+k}\\
x_{M} & x_{M-1}%
\end{array}
\right)  }\cdot\underbrace{B}_{=\left(
\begin{array}
[c]{cc}%
a & 1\\
b & 0
\end{array}
\right)  }\nonumber\\
&  =\left(
\begin{array}
[c]{cc}%
x_{M+k+1} & x_{M+k}\\
x_{M} & x_{M-1}%
\end{array}
\right)  \cdot\left(
\begin{array}
[c]{cc}%
a & 1\\
b & 0
\end{array}
\right) \nonumber\\
&  =\left(
\begin{array}
[c]{cc}%
x_{M+k+1}\cdot a+x_{M+k}\cdot b & x_{M+k+1}\cdot1+x_{M+k}\cdot0\\
x_{M}\cdot a+x_{M-1}\cdot b & x_{M}\cdot1+x_{M-1}\cdot0
\end{array}
\right) \nonumber\\
&  \ \ \ \ \ \ \ \ \ \ \left(  \text{by the definition of a product of two
matrices}\right) \nonumber\\
&  =\left(
\begin{array}
[c]{cc}%
ax_{M+k+1}+bx_{M+k} & x_{M+k+1}\\
ax_{M}+bx_{M-1} & x_{M}%
\end{array}
\right)  . \label{sol.ps4.det.fibo.BmC.4}%
\end{align}

But (\ref{eq.det.fibo.rec}) (applied to $n=M+k+2$) yields%
\[
x_{M+k+2}=a\underbrace{x_{\left(  M+k+2\right)  -1}}_{=x_{M+k+1}%
}+b\underbrace{x_{\left(  M+k+2\right)  -2}}_{=x_{M+k}}=ax_{M+k+1}+bx_{M+k}.
\]
Also, (\ref{eq.det.fibo.rec}) (applied to $n=M+1$) yields $x_{M+1}%
=a\underbrace{x_{\left(  M+1\right)  -1}}_{=x_{M}}+b\underbrace{x_{\left(
M+1\right)  -2}}_{=x_{M-1}}=ax_{M}+bx_{M-1}$. Now,%
\[
\left(
\begin{array}
[c]{cc}%
x_{M+k+2} & x_{M+k+1}\\
x_{M+1} & x_{M}%
\end{array}
\right)  =\left(
\begin{array}
[c]{cc}%
ax_{M+k+1}+bx_{M+k} & x_{M+k+1}\\
ax_{M}+bx_{M-1} & x_{M}%
\end{array}
\right)
\]
(since $x_{M+k+2}=ax_{M+k+1}+bx_{M+k}$ and $x_{M+1}=ax_{M}+bx_{M-1}$).
Compared with (\ref{sol.ps4.det.fibo.BmC.4}), this yields $AB^{M}=\left(
\begin{array}
[c]{cc}%
x_{M+k+2} & x_{M+k+1}\\
x_{M+1} & x_{M}%
\end{array}
\right)  $. In other words, (\ref{sol.ps4.det.fibo.BmC}) holds for $m=M$. This
completes the induction step. Thus, (\ref{sol.ps4.det.fibo.BmC}) is proven by induction.]

Now, let $n>k$ be an integer. Then, $n-k>0$, so that $n-k-1\in\mathbb{N}$.
Hence, (\ref{sol.ps4.det.fibo.BmC}) (applied to $m=n-k-1$) yields%
\[
AB^{n-k-1}=\left(
\begin{array}
[c]{cc}%
x_{\left(  n-k-1\right)  +k+2} & x_{\left(  n-k-1\right)  +k+1}\\
x_{\left(  n-k-1\right)  +1} & x_{n-k-1}%
\end{array}
\right)  =\left(
\begin{array}
[c]{cc}%
x_{n+1} & x_{n}\\
x_{n-k} & x_{n-k-1}%
\end{array}
\right)  .
\]
Taking determinants on both sides of this equality, we obtain%
\[
\det\left(  AB^{n-k-1}\right)  =\det\left(
\begin{array}
[c]{cc}%
x_{n+1} & x_{n}\\
x_{n-k} & x_{n-k-1}%
\end{array}
\right)  =x_{n+1}x_{n-k-1}-x_{n}x_{n-k}.
\]
Hence,%
\begin{align*}
&  x_{n+1}x_{n-k-1}-x_{n}x_{n-k}\\
&  =\det\left(  AB^{n-k-1}\right)  =\det A\cdot\underbrace{\det\left(
B^{n-k-1}\right)  }_{\substack{=\left(  \det B\right)  ^{n-k-1}\\\text{(by
Corollary \ref{cor.det.product} \textbf{(b)}, applied}\\\text{to }2\text{ and
}n-k-1\text{ instead of }n\text{ and }k\text{)}}}\\
&  \ \ \ \ \ \ \ \ \ \ \left(  \text{by Theorem \ref{thm.det(AB)}, applied to
}2\text{ and }B^{n-k-1}\text{ instead of }n\text{ and }B\right) \\
&  =\underbrace{\det A}_{=x_{k+2}x_{0}-x_{k+1}x_{1}}\cdot\left(
\underbrace{\det B}_{=-b}\right)  ^{n-k-1}=\left(  x_{k+2}x_{0}-x_{k+1}%
x_{1}\right)  \cdot\left(  -b\right)  ^{n-k-1}\\
&  =\left(  -b\right)  ^{n-k-1}\left(  x_{k+2}x_{0}-x_{k+1}x_{1}\right)  .
\end{align*}
This solves Exercise \ref{exe.ps4.det.fibo}.
\end{proof}

\subsection{Solution to Exercise \ref{exe.ps4.pascal}}

\begin{proof}
[Solution to Exercise \ref{exe.ps4.pascal}.]We let $B$ be the $n\times
n$-matrix $\left(  \dbinom{i-1}{j-1}\right)  _{1\leq i\leq n,\ 1\leq j\leq n}%
$. We have $\dbinom{i-1}{j-1}=0$ for every $\left(  i,j\right)  \in\left\{
1,2,\ldots,n\right\}  ^{2}$ satisfying $i<j$\ \ \ \ \footnote{\textit{Proof.}
Let $\left(  i,j\right)  \in\left\{  1,2,\ldots,n\right\}  ^{2}$ be such that
$i<j$. Then, $i-1\in\mathbb{N}$ (since $i\geq1$) and $j-1\in\mathbb{N}$ (since
$j\geq1$) and $i-1<j-1$ (since $i<j$). Hence, (\ref{eq.binom.0}) (applied to
$i-1$ and $j-1$ instead of $m$ and $n$) yields $\dbinom{i-1}{j-1}=0$, qed.}.
Therefore, Exercise \ref{exe.ps4.3} (applied to $B$ and $\dbinom{i-1}{j-1}$
instead of $A$ and $a_{i,j}$) yields
\[
\det B=\underbrace{\dbinom{1-1}{1-1}}_{=1}\underbrace{\dbinom{2-1}{2-1}}%
_{=1}\cdots\underbrace{\dbinom{n-1}{n-1}}_{=1}=1\cdot1\cdot\cdots\cdot1=1.
\]
Exercise \ref{exe.ps4.4} (applied to $B$ instead of $A$) shows that
$\det\left(  B^{T}\right)  =\det B=1$.

Now, we are going to show that $A=BB^{T}$.

Indeed, let us show that%
\begin{equation}
\sum_{k=1}^{n}\dbinom{i-1}{k-1}\dbinom{j-1}{k-1}=\dbinom{i+j-2}{i-1}
\label{sol.ps4.pascal.1}%
\end{equation}
for every $\left(  i,j\right)  \in\left\{  1,2,\ldots,n\right\}  ^{2}$.

[\textit{Proof of (\ref{sol.ps4.pascal.1}):} Let $\left(  i,j\right)
\in\left\{  1,2,\ldots,n\right\}  ^{2}$. Then, $i\geq1$, so that
$i-1\in\mathbb{N}$. Hence, Proposition \ref{prop.vandermonde.consequences}
\textbf{(b)} (applied to $x=i-1$ and $y=j-1$) yields
\begin{equation}
\dbinom{\left(  i-1\right)  +\left(  j-1\right)  }{i-1}=\sum_{k=0}%
^{i-1}\dbinom{i-1}{k}\dbinom{j-1}{k}. \label{sol.ps4.pascal.2}%
\end{equation}
On the other hand, for every $k\in\left\{  i+1,i+2,\ldots,n\right\}  $, we
have%
\begin{equation}
\dbinom{i-1}{k-1}=0 \label{sol.ps4.pascal.4}%
\end{equation}
\footnote{\textit{Proof of (\ref{sol.ps4.pascal.4}):} Let $k\in\left\{
i+1,i+2,\ldots,n\right\}  $. Then, $k>i\geq1$ and thus $k-1\in\mathbb{N}$.
Also, $k>i$, so that $k-i>i-1$ and thus $i-1<k-1$. Hence, (\ref{eq.binom.0})
(applied to $i-1$ and $k-1$ instead of $m$ and $n$) yields $\dbinom{i-1}%
{k-1}=0$, qed.}. Now, $i\leq n$, so that%
\begin{align*}
&  \sum_{k=1}^{n}\dbinom{i-1}{k-1}\dbinom{j-1}{k-1}\\
&  =\sum_{k=1}^{i}\dbinom{i-1}{k-1}\dbinom{j-1}{k-1}+\sum_{k=i+1}%
^{n}\underbrace{\dbinom{i-1}{k-1}}_{\substack{=0\\\text{(by
(\ref{sol.ps4.pascal.4}))}}}\dbinom{j-1}{k-1}\\
&  =\sum_{k=1}^{i}\dbinom{i-1}{k-1}\dbinom{j-1}{k-1}+\underbrace{\sum
_{k=i+1}^{n}0\dbinom{j-1}{k-1}}_{=0}=\sum_{k=1}^{i}\dbinom{i-1}{k-1}%
\dbinom{j-1}{k-1}\\
&  =\sum_{k=0}^{i-1}\dbinom{i-1}{k}\dbinom{j-1}{k}\\
&  \ \ \ \ \ \ \ \ \ \ \left(  \text{here, we have substituted }k\text{ for
}k-1\text{ in the sum}\right) \\
&  =\dbinom{\left(  i-1\right)  +\left(  j-1\right)  }{i-1}%
\ \ \ \ \ \ \ \ \ \ \left(  \text{by (\ref{sol.ps4.pascal.2})}\right) \\
&  =\dbinom{i+j-2}{i-1}\ \ \ \ \ \ \ \ \ \ \left(  \text{since }\left(
i-1\right)  +\left(  j-1\right)  =i+j-2\right)  .
\end{align*}
This proves (\ref{sol.ps4.pascal.1}).]

Now, we have $B=\left(  \dbinom{i-1}{j-1}\right)  _{1\leq i\leq n,\ 1\leq
j\leq n}$ and therefore \newline$B^{T}=\left(  \dbinom{j-1}{i-1}\right)
_{1\leq i\leq n,\ 1\leq j\leq n}$ (by the definition of $B^{T}$). Hence, the
definition of the product $BB^{T}$ shows that%
\[
BB^{T}=\left(  \underbrace{\sum_{k=1}^{n}\dbinom{i-1}{k-1}\dbinom{j-1}{k-1}%
}_{\substack{=\dbinom{i+j-2}{i-1}\\\text{(by (\ref{sol.ps4.pascal.1}))}%
}}\right)  _{1\leq i\leq n,\ 1\leq j\leq n}=\left(  \dbinom{i+j-2}%
{i-1}\right)  _{1\leq i\leq n,\ 1\leq j\leq n}=A.
\]
Hence, $A=BB^{T}$, so that%
\begin{align*}
\det A  &  =\det\left(  BB^{T}\right)  =\underbrace{\det B}_{=1}%
\cdot\underbrace{\det\left(  B^{T}\right)  }_{=1}\\
&  \ \ \ \ \ \ \ \ \ \ \left(  \text{by Theorem \ref{thm.det(AB)}, applied to
}B\text{ and }B^{T}\text{ instead of }A\text{ and }B\right) \\
&  =1.
\end{align*}
This solves Exercise \ref{exe.ps4.pascal}.
\end{proof}

\subsection{Solution to Exercise \ref{exe.lem.increasing-sequences}}

\begin{proof}
[Proof of Lemma \ref{lem.increasing-sequences}.]\textbf{(b)} Let $\left(
g_{1},g_{2},\ldots,g_{n}\right)  \in\left\{  1,2,\ldots,m\right\}  ^{n}$ be an
$n$-tuple satisfying $g_{1}<g_{2}<\cdots<g_{n}$. We shall derive a contradiction.

We have $m<n$, so that $n>m\geq0$. Thus, $n\geq1$ (since $n$ is an integer).
Therefore, the elements $g_{1}$ and $g_{n}$ of $\left\{  1,2,\ldots,m\right\}
$ are well-defined. Every $i\in\left\{  1,2,\ldots,n-1\right\}  $ satisfies
$g_{i}-i\leq g_{i+1}-\left(  i+1\right)  $\ \ \ \ \footnote{\textit{Proof.}
Let $i\in\left\{  1,2,\ldots,n-1\right\}  $. Then, $g_{i}<g_{i+1}$ (since
$g_{1}<g_{2}<\cdots<g_{n}$). Hence, $g_{i}\leq g_{i+1}-1$ (since $g_{i}$ and
$g_{i+1}$ are integers). Hence, $\underbrace{g_{i}}_{\leq g_{i+1}-1}-i\leq
g_{i+1}-1-i=g_{i+1}-\left(  i+1\right)  $, qed.}. In other words, we have
$g_{1}-1\leq g_{2}-2\leq\cdots\leq g_{n}-n$. Hence, $g_{1}-1\leq g_{n}-n$.

But $g_{n}\in\left\{  1,2,\ldots,m\right\}  $, so that $g_{n}\leq m<n$ and
thus $\underbrace{g_{n}}_{<n}-n<n-n=0$. Hence, $g_{1}-1\leq g_{n}-n<0$. On the
other hand, $g_{1}\in\left\{  1,2,\ldots,m\right\}  $, so that $g_{1}\geq1$
and thus $g_{1}-1\geq0$. This contradicts $g_{1}-1<0$.

Now, let us forget that we fixed $\left(  g_{1},g_{2},\ldots,g_{n}\right)  $.
We thus have derived a contradiction for every $n$-tuple $\left(  g_{1}%
,g_{2},\ldots,g_{n}\right)  \in\left\{  1,2,\ldots,m\right\}  ^{n}$ satisfying
$g_{1}<g_{2}<\cdots<g_{n}$. Hence, there exists no $n$-tuple $\left(
g_{1},g_{2},\ldots,g_{n}\right)  \in\left\{  1,2,\ldots,m\right\}  ^{n}$
satisfying $g_{1}<g_{2}<\cdots<g_{n}$. This proves Lemma
\ref{lem.increasing-sequences} \textbf{(b)}.

\textbf{(a)} The $n$-tuple $\left(  1,2,\ldots,n\right)  $ clearly belongs to
$\left\{  1,2,\ldots,n\right\}  ^{n}$, and satisfies $1<2<\cdots<n$. Hence,
there exists an $n$-tuple $\left(  g_{1},g_{2},\ldots,g_{n}\right)
\in\left\{  1,2,\ldots,n\right\}  ^{n}$ satisfying $g_{1}<g_{2}<\cdots<g_{n}$,
namely the $n$-tuple $\left(  1,2,\ldots,n\right)  $. We shall now show that
$\left(  1,2,\ldots,n\right)  $ is the only such $n$-tuple.

Indeed, let $\left(  g_{1},g_{2},\ldots,g_{n}\right)  \in\left\{
1,2,\ldots,n\right\}  ^{n}$ be an $n$-tuple satisfying $g_{1}<g_{2}%
<\cdots<g_{n}$. We shall prove that $\left(  g_{1},g_{2},\ldots,g_{n}\right)
=\left(  1,2,\ldots,n\right)  $.

Every $i\in\left\{  1,2,\ldots,n-1\right\}  $ satisfies $g_{i}-i\leq
g_{i+1}-\left(  i+1\right)  $\ \ \ \ \footnote{\textit{Proof.} Let
$i\in\left\{  1,2,\ldots,n-1\right\}  $. Then, $g_{i}<g_{i+1}$ (since
$g_{1}<g_{2}<\cdots<g_{n}$). Hence, $g_{i}\leq g_{i+1}-1$ (since $g_{i}$ and
$g_{i+1}$ are integers). Hence, $\underbrace{g_{i}}_{\leq g_{i+1}-1}-i\leq
g_{i+1}-1-i=g_{i+1}-\left(  i+1\right)  $, qed.}. In other words, we have
$g_{1}-1\leq g_{2}-2\leq\cdots\leq g_{n}-n$. In other words,%
\begin{equation}
g_{u}-u\leq g_{v}-v \label{pf.lem.increasing-sequences.b.1}%
\end{equation}
for any two elements $u$ and $v$ of $\left\{  1,2,\ldots,n\right\}  $
satisfying $u\leq v$.

Let $i\in\left\{  1,2,\ldots,n\right\}  $. Then, $1\leq i\leq n$, so that
$1\leq n$. Hence, the elements $g_{1}$ and $g_{n}$ of $\left\{  1,2,\ldots
,n\right\}  $ are well-defined.

Now, $1\leq i$. Hence, (\ref{pf.lem.increasing-sequences.b.1}) (applied to
$u=1$ and $v=i$) yields $g_{1}-1\leq g_{i}-i$. But $g_{1}\in\left\{
1,2,\ldots,n\right\}  $, so that $g_{1}\geq1$ and thus $g_{1}-1\geq0$. Hence,
$0\leq g_{1}-1\leq g_{i}-i$.

On the other hand, $i\leq n$. Therefore,
(\ref{pf.lem.increasing-sequences.b.1}) (applied to $u=i$ and $v=n$) yields
$g_{i}-i\leq g_{n}-n$. But $g_{n}\in\left\{  1,2,\ldots,n\right\}  $, so that
$g_{n}\leq n$ and thus $g_{n}-n\leq0$. Hence, $g_{i}-i\leq g_{n}-n\leq0$.
Combined with $0\leq g_{i}-i$, this yields $g_{i}-i=0$, so that $g_{i}=i$.

Now, let us forget that we fixed $i$. We thus have shown that $g_{i}=i$ for
every $i\in\left\{  1,2,\ldots,n\right\}  $. In other words, $\left(
g_{1},g_{2},\ldots,g_{n}\right)  =\left(  1,2,\ldots,n\right)  $.

Let us now forget that we fixed $\left(  g_{1},g_{2},\ldots,g_{n}\right)  $.
We thus have shown that if $\left(  g_{1},g_{2},\ldots,g_{n}\right)
\in\left\{  1,2,\ldots,n\right\}  ^{n}$ is an $n$-tuple satisfying
$g_{1}<g_{2}<\cdots<g_{n}$, then $\left(  g_{1},g_{2},\ldots,g_{n}\right)
=\left(  1,2,\ldots,n\right)  $. In other words, every $n$-tuple $\left(
g_{1},g_{2},\ldots,g_{n}\right)  \in\left\{  1,2,\ldots,n\right\}  ^{n}$
satisfying $g_{1}<g_{2}<\cdots<g_{n}$ must be equal to $\left(  1,2,\ldots
,n\right)  $. Hence, there exists at most one $n$-tuple $\left(  g_{1}%
,g_{2},\ldots,g_{n}\right)  \in\left\{  1,2,\ldots,n\right\}  ^{n}$ satisfying
$g_{1}<g_{2}<\cdots<g_{n}$ (namely, the $n$-tuple $\left(  1,2,\ldots
,n\right)  $).

We now know the following two facts:

\begin{itemize}
\item There exists an $n$-tuple $\left(  g_{1},g_{2},\ldots,g_{n}\right)
\in\left\{  1,2,\ldots,n\right\}  ^{n}$ satisfying $g_{1}<g_{2}<\cdots<g_{n}$,
namely the $n$-tuple $\left(  1,2,\ldots,n\right)  $.

\item There exists at most one $n$-tuple $\left(  g_{1},g_{2},\ldots
,g_{n}\right)  \in\left\{  1,2,\ldots,n\right\}  ^{n}$ satisfying $g_{1}%
<g_{2}<\cdots<g_{n}$.
\end{itemize}

Combining these two facts, we conclude that there exists exactly one $n$-tuple
$\left(  g_{1},g_{2},\ldots,g_{n}\right)  \in\left\{  1,2,\ldots,n\right\}
^{n}$ satisfying $g_{1}<g_{2}<\cdots<g_{n}$, namely the $n$-tuple $\left(
1,2,\ldots,n\right)  $. This proves Lemma \ref{lem.increasing-sequences}
\textbf{(a)}.
\end{proof}

\subsection{Solution to Exercise \ref{exe.sorting.basics}}

Before we give a formal proof of Proposition \ref{prop.sorting}, let us
outline the main ideas of this proof: We shall define the notion of
\textit{inversion} of an $n$-tuple $\left(  a_{1},a_{2},\ldots,a_{n}\right)  $
of integers\footnote{It will be defined as a pair $\left(  i,j\right)  $ of
integers satisfying $1\leq i<j\leq n$ and $a_{i}>a_{j}$. The analogy with the
notion of \textquotedblleft inversion\textquotedblright\ of a permutation is
intentional.}; then we will argue that swapping two adjacent entries $a_{k}$
and $a_{k+1}$ of an $n$-tuple $\left(  a_{1},a_{2},\ldots,a_{n}\right)  $
which are \textquotedblleft out of order\textquotedblright\ (i.e., satisfy
$a_{k}>a_{k+1}$) reduces the number of inversions by $1$ (this is our equality
(\ref{pf.prop.sorting.a.fact5}) further below), and thus a sequence of such
swaps will eventually end and therefore bring the tuple into weakly increasing
order. (This is similar to an argument in the solution of Exercise
\ref{exe.ps2.2.5} \textbf{(e)}.) Other proofs of Proposition
\ref{prop.sorting} \textbf{(a)} are possible\footnote{Roughly speaking, to
each \href{https://en.wikipedia.org/wiki/Sorting_algorithm}{sorting algorithm}
corresponds at least one proof of Proposition \ref{prop.sorting} \textbf{(a)}.
(\textquotedblleft At least\textquotedblright\ because there are often several
ways to prove the correctness of a given sorting algorithm.) The proof we just
outlined corresponds to \textquotedblleft bubble sort\textquotedblright.}.
Proposition \ref{prop.sorting} \textbf{(b)} will then easily follow (indeed,
we will argue that $a_{\sigma\left(  i\right)  }$ is the smallest integer $x$
such that at least $i$ different elements $j\in\left\{  1,2,\ldots,n\right\}
$ satisfy $a_{j}\leq x$), and Proposition \ref{prop.sorting} \textbf{(c)} will
finally follow from parts \textbf{(a)} and \textbf{(b)}.

Here is the proof in detail:

\begin{proof}
[Proof of Proposition \ref{prop.sorting}.]\textbf{(a)} Let us forget that we
fixed $a_{1},a_{2},\ldots,a_{n}$.

We shall first introduce some notations.

We let $\left[  n\right]  $ denote the set $\left\{  1,2,\ldots,n\right\}  $.

If $\mathbf{a}=\left(  a_{1},a_{2},\ldots,a_{n}\right)  $ is an $n$-tuple of
integers, then:

\begin{itemize}
\item An \textit{inversion} of $\mathbf{a}$ will mean a pair $\left(
i,j\right)  \in\left[  n\right]  ^{2}$ satisfying $i<j$ and $a_{i}>a_{j}$.

\item We denote by $\operatorname*{Inv}\left(  \mathbf{a}\right)  $ the set of
all inversions of $\mathbf{a}$. Thus, $\operatorname*{Inv}\left(
\mathbf{a}\right)  \subseteq\left[  n\right]  ^{2}$. More precisely,%
\begin{align}
\operatorname*{Inv}\left(  \mathbf{a}\right)   &  =\left(  \text{the set of
all inversions of }\mathbf{a}\right) \nonumber\\
&  =\left\{  \left(  i,j\right)  \in\left[  n\right]  ^{2}\ \mid\ i<j\text{
and }a_{i}>a_{j}\right\} \label{pf.prop.sorting.a.Inv.1}\\
&  \ \ \ \ \ \ \ \ \ \ \left(
\begin{array}
[c]{c}%
\text{since the inversions of }\mathbf{a}\text{ are the pairs }\left(
i,j\right)  \in\left[  n\right]  ^{2}\text{ satisfying}\\
i<j\text{ and }a_{i}>a_{j}\text{ (by the definition of an \textquotedblleft
inversion\textquotedblright)}%
\end{array}
\right) \nonumber\\
&  =\left\{  \left(  u,v\right)  \in\left[  n\right]  ^{2}\ \mid\ u<v\text{
and }a_{u}>a_{v}\right\} \label{pf.prop.sorting.a.Inv.2}\\
&  \ \ \ \ \ \ \ \ \ \ \left(  \text{here, we renamed the index }\left(
i,j\right)  \text{ as }\left(  u,v\right)  \right)  .\nonumber
\end{align}

\item We denote by $\ell\left(  \mathbf{a}\right)  $ the number $\left\vert
\operatorname*{Inv}\left(  \mathbf{a}\right)  \right\vert $. (This is
well-defined because $\operatorname*{Inv}\left(  \mathbf{a}\right)  $ is
finite (since $\operatorname*{Inv}\left(  \mathbf{a}\right)  \subseteq\left[
n\right]  ^{2}$).) Thus,%
\begin{align}
\ell\left(  \mathbf{a}\right)   &  =\left\vert \underbrace{\operatorname*{Inv}%
\left(  \mathbf{a}\right)  }_{=\left(  \text{the set of all inversions of
}\mathbf{a}\right)  }\right\vert =\left\vert \left(  \text{the set of all
inversions of }\mathbf{a}\right)  \right\vert \nonumber\\
&  =\left(  \text{the number of all inversions of }\mathbf{a}\right)  .
\label{pf.prop.sorting.a.l(sigma)}%
\end{align}

\item For every permutation $\tau\in S_{n}$, we denote by $\mathbf{a}\circ
\tau$ the $n$-tuple \newline$\left(  a_{\tau\left(  1\right)  },a_{\tau\left(
2\right)  },\ldots,a_{\tau\left(  n\right)  }\right)  $ of integers.
\end{itemize}

We notice that if $\mathbf{a}$ is any $n$-tuple of integers, then%
\begin{equation}
\mathbf{a}\circ\left(  \sigma\circ\tau\right)  =\left(  \mathbf{a}\circ
\sigma\right)  \circ\tau\ \ \ \ \ \ \ \ \ \ \text{for any }\sigma\in
S_{n}\text{ and }\tau\in S_{n} \label{pf.prop.sorting.a.sigmatau}%
\end{equation}
\footnote{\textit{Proof of (\ref{pf.prop.sorting.a.sigmatau}):} Let
$\mathbf{a}$ be any $n$-tuple of integers. Let $\sigma\in S_{n}$ and $\tau\in
S_{n}$. We must show that $\mathbf{a}\circ\left(  \sigma\circ\tau\right)
=\left(  \mathbf{a}\circ\sigma\right)  \circ\tau$.
\par
Write the $n$-tuple $\mathbf{a}$ in the form $\mathbf{a}=\left(  a_{1}%
,a_{2},\ldots,a_{n}\right)  $ for some integers $a_{1},a_{2},\ldots,a_{n}$.
Thus, the definition of $\mathbf{a}\circ\left(  \sigma\circ\tau\right)  $
yields%
\[
\mathbf{a}\circ\left(  \sigma\circ\tau\right)  =\left(  a_{\left(  \sigma
\circ\tau\right)  \left(  1\right)  },a_{\left(  \sigma\circ\tau\right)
\left(  2\right)  },\ldots,a_{\left(  \sigma\circ\tau\right)  \left(
n\right)  }\right)  =\left(  a_{\sigma\left(  \tau\left(  1\right)  \right)
},a_{\sigma\left(  \tau\left(  2\right)  \right)  },\ldots,a_{\sigma\left(
\tau\left(  n\right)  \right)  }\right)
\]
(since $a_{\left(  \sigma\circ\tau\right)  \left(  i\right)  }=a_{\sigma
\left(  \tau\left(  i\right)  \right)  }$ for every $i\in\left\{
1,2,\ldots,n\right\}  $). On the other hand, the definition of $\mathbf{a}%
\circ\sigma$ yields $\mathbf{a}\circ\sigma=\left(  a_{\sigma\left(  1\right)
},a_{\sigma\left(  2\right)  },\ldots,a_{\sigma\left(  n\right)  }\right)  $
(since $\mathbf{a}=\left(  a_{1},a_{2},\ldots,a_{n}\right)  $). Therefore, the
definition of $\left(  \mathbf{a}\circ\sigma\right)  \circ\tau$ yields
$\left(  \mathbf{a}\circ\sigma\right)  \circ\tau=\left(  a_{\sigma\left(
\tau\left(  1\right)  \right)  },a_{\sigma\left(  \tau\left(  2\right)
\right)  },\ldots,a_{\sigma\left(  \tau\left(  n\right)  \right)  }\right)  $.
Compared with $\mathbf{a}\circ\left(  \sigma\circ\tau\right)  =\left(
a_{\sigma\left(  \tau\left(  1\right)  \right)  },a_{\sigma\left(  \tau\left(
2\right)  \right)  },\ldots,a_{\sigma\left(  \tau\left(  n\right)  \right)
}\right)  $, this yields $\mathbf{a}\circ\left(  \sigma\circ\tau\right)
=\left(  \mathbf{a}\circ\sigma\right)  \circ\tau$. This proves
(\ref{pf.prop.sorting.a.sigmatau}).}. Also, if $\mathbf{a}$ is any $n$-tuple
of integers, then%
\begin{equation}
\mathbf{a}\circ\operatorname*{id}=\mathbf{a} \label{pf.prop.sorting.a.id}%
\end{equation}
\footnote{\textit{Proof of (\ref{pf.prop.sorting.a.id}):} Let $\mathbf{a}$ be
any $n$-tuple of integers.
\par
Write the $n$-tuple $\mathbf{a}$ in the form $\mathbf{a}=\left(  a_{1}%
,a_{2},\ldots,a_{n}\right)  $ for some integers $a_{1},a_{2},\ldots,a_{n}$.
Thus, the definition of $\mathbf{a}\circ\operatorname*{id}$ yields%
\[
\mathbf{a}\circ\operatorname*{id}=\left(  a_{\operatorname*{id}\left(
1\right)  },a_{\operatorname*{id}\left(  2\right)  },\ldots
,a_{\operatorname*{id}\left(  n\right)  }\right)  =\left(  a_{1},a_{2}%
,\ldots,a_{n}\right)
\]
(since $a_{\operatorname*{id}\left(  i\right)  }=a_{i}$ for every
$i\in\left\{  1,2,\ldots,n\right\}  $). Compared with $\mathbf{a}=\left(
a_{1},a_{2},\ldots,a_{n}\right)  $, this yields $\mathbf{a}\circ
\operatorname*{id}=\mathbf{a}$. This proves (\ref{pf.prop.sorting.a.id}).}.

If $X$, $X^{\prime}$, $Y$ and $Y^{\prime}$ are four sets and if $\alpha
:X\rightarrow X^{\prime}$ and $\beta:Y\rightarrow Y^{\prime}$ are two maps,
then $\alpha\times\beta$ will denote the map%
\begin{align*}
X\times Y  &  \rightarrow X^{\prime}\times Y^{\prime},\\
\left(  x,y\right)   &  \mapsto\left(  \alpha\left(  x\right)  ,\beta\left(
y\right)  \right)  .
\end{align*}

Recall that, for each $k\in\left\{  1,2,\ldots,n-1\right\}  $, we have defined
$s_{k}$ to be the permutation in $S_{n}$ that swaps $k$ with $k+1$ but leaves
all other numbers unchanged. This permutation $s_{k}$ is a map $\left[
n\right]  \rightarrow\left[  n\right]  $ and satisfies $s_{k}^{2}%
=\operatorname*{id}$. For every $k\in\left\{  1,2,\ldots,n-1\right\}  $, the
map $s_{k}\times s_{k}:\left[  n\right]  \times\left[  n\right]
\rightarrow\left[  n\right]  \times\left[  n\right]  $ satisfies $\left(
s_{k}\times s_{k}\right)  ^{2}=\operatorname*{id}$%
\ \ \ \ \footnote{\textit{Proof.} Let $u\in\left[  n\right]  \times\left[
n\right]  $. Then, we can write $u$ in the form $\left(  i,j\right)  $ for
some $i\in\left[  n\right]  $ and $j\in\left[  n\right]  $. Consider these $i$
and $j$. We have%
\begin{align*}
\underbrace{\left(  s_{k}\times s_{k}\right)  ^{2}}_{=\left(  s_{k}\times
s_{k}\right)  \circ\left(  s_{k}\times s_{k}\right)  }\left(  \underbrace{u}%
_{=\left(  i,j\right)  }\right)   &  =\left(  \left(  s_{k}\times
s_{k}\right)  \circ\left(  s_{k}\times s_{k}\right)  \right)  \left(  \left(
i,j\right)  \right)  =\left(  s_{k}\times s_{k}\right)  \left(
\underbrace{\left(  s_{k}\times s_{k}\right)  \left(  \left(  i,j\right)
\right)  }_{\substack{=\left(  s_{k}\left(  i\right)  ,s_{k}\left(  j\right)
\right)  \\\text{(by the definition of }s_{k}\times s_{k}\text{)}}}\right) \\
&  =\left(  s_{k}\times s_{k}\right)  \left(  \left(  s_{k}\left(  i\right)
,s_{k}\left(  j\right)  \right)  \right)  =\left(  \underbrace{s_{k}\left(
s_{k}\left(  i\right)  \right)  }_{=s_{k}^{2}\left(  i\right)  }%
,\underbrace{s_{k}\left(  s_{k}\left(  j\right)  \right)  }_{=s_{k}^{2}\left(
j\right)  }\right) \\
&  \ \ \ \ \ \ \ \ \ \ \left(  \text{by the definition of }s_{k}\times
s_{k}\right) \\
&  =\left(  \underbrace{s_{k}^{2}}_{=\operatorname*{id}}\left(  i\right)
,\underbrace{s_{k}^{2}}_{=\operatorname*{id}}\left(  j\right)  \right)
=\left(  \underbrace{\operatorname*{id}\left(  i\right)  }_{=i}%
,\underbrace{\operatorname*{id}\left(  j\right)  }_{=j}\right)  =\left(
i,j\right)  =u=\operatorname*{id}\left(  u\right)  .
\end{align*}
\par
Now, let us forget that we fixed $u$. We thus have shown that $\left(
s_{k}\times s_{k}\right)  ^{2}\left(  u\right)  =\operatorname*{id}\left(
u\right)  $ for every $u\in\left[  n\right]  \times\left[  n\right]  $. In
other words, $\left(  s_{k}\times s_{k}\right)  ^{2}=\operatorname*{id}$,
qed.}, and thus is a bijection from $\left[  n\right]  ^{2}$ to $\left[
n\right]  ^{2}$\ \ \ \ \footnote{\textit{Proof.} Let $k\in\left\{
1,2,\ldots,n-1\right\}  $. We have $\left(  s_{k}\times s_{k}\right)
\circ\left(  s_{k}\times s_{k}\right)  =\left(  s_{k}\times s_{k}\right)
^{2}=\operatorname*{id}$. Hence, the maps $s_{k}\times s_{k}$ and $s_{k}\times
s_{k}$ are mutually inverse. Hence, the map $s_{k}\times s_{k}$ is invertible,
thus a bijection. Therefore, this map $s_{k}\times s_{k}$ is a bijection from
$\left[  n\right]  ^{2}$ to $\left[  n\right]  ^{2}$ (because its domain is
$\left[  n\right]  \times\left[  n\right]  =\left[  n\right]  ^{2}$, and its
codomain is $\left[  n\right]  \times\left[  n\right]  =\left[  n\right]
^{2}$). Qed.}.

Let us now recall a simple fact: If $u$ and $v$ are two integers such that
$1\leq u<v\leq n$, and if $k\in\left\{  1,2,\ldots,n-1\right\}  $ is such that
$\left(  u,v\right)  \neq\left(  k,k+1\right)  $, then%
\begin{equation}
s_{k}\left(  u\right)  <s_{k}\left(  v\right)  . \label{pf.prop.sorting.a.2.5}%
\end{equation}
(This was proven in the solution of Exercise \ref{exe.ps2.2.5} \textbf{(a)}.)

Now, we notice the following facts:

\begin{itemize}
\item If $\mathbf{a}=\left(  a_{1},a_{2},\ldots,a_{n}\right)  $ is an
$n$-tuple of integers, and if $k\in\left\{  1,2,\ldots,n-1\right\}  $, then%
\begin{equation}
\operatorname*{Inv}\left(  \mathbf{a}\circ s_{k}\right)  =\left\{  \left(
u,v\right)  \in\left[  n\right]  ^{2}\ \mid\ u<v\text{ and }a_{s_{k}\left(
u\right)  }>a_{s_{k}\left(  v\right)  }\right\}
\label{pf.prop.sorting.a.fact1}%
\end{equation}
\footnote{\textit{Proof of (\ref{pf.prop.sorting.a.fact1}):} Let
$\mathbf{a}=\left(  a_{1},a_{2},\ldots,a_{n}\right)  $ be an $n$-tuple of
integers, and let $k\in\left\{  1,2,\ldots,n-1\right\}  $. The definition of
$\mathbf{a}\circ s_{k}$ yields $\mathbf{a}\circ s_{k}=\left(  a_{s_{k}\left(
1\right)  },a_{s_{k}\left(  2\right)  },\ldots,a_{s_{k}\left(  n\right)
}\right)  $ (since $\mathbf{a}=\left(  a_{1},a_{2},\ldots,a_{n}\right)  $).
Hence, (\ref{pf.prop.sorting.a.Inv.2}) (applied to $\mathbf{a}\circ s_{k}$ and
$a_{s_{k}\left(  i\right)  }$ instead of $\mathbf{a}$ and $a_{i}$) yields
$\operatorname*{Inv}\left(  \mathbf{a}\circ s_{k}\right)  =\left\{  \left(
u,v\right)  \in\left[  n\right]  ^{2}\ \mid\ u<v\text{ and }a_{s_{k}\left(
u\right)  }>a_{s_{k}\left(  v\right)  }\right\}  $. This proves
(\ref{pf.prop.sorting.a.fact1}).}.

\item If $\mathbf{a}$ is an $n$-tuple of integers, and if $k\in\left\{
1,2,\ldots,n-1\right\}  $, then%
\begin{equation}
\left(  s_{k}\times s_{k}\right)  ^{-1}\left(  \operatorname*{Inv}\left(
\mathbf{a}\right)  \setminus\left\{  \left(  k,k+1\right)  \right\}  \right)
\subseteq\operatorname*{Inv}\left(  \mathbf{a}\circ s_{k}\right)
\setminus\left\{  \left(  k,k+1\right)  \right\}
\label{pf.prop.sorting.a.fact2}%
\end{equation}
\footnote{\textit{Proof of (\ref{pf.prop.sorting.a.fact2}):} Let $\mathbf{a}$
be an $n$-tuple of integers, and let $k\in\left\{  1,2,\ldots,n-1\right\}  $.
Write the $n$-tuple $\mathbf{a}$ in the form $\mathbf{a}=\left(  a_{1}%
,a_{2},\ldots,a_{n}\right)  $ for some integers $a_{1},a_{2},\ldots,a_{n}$.
\par
Let $c\in\left(  s_{k}\times s_{k}\right)  ^{-1}\left(  \operatorname*{Inv}%
\left(  \mathbf{a}\right)  \setminus\left\{  \left(  k,k+1\right)  \right\}
\right)  $. Thus, $c\in\left[  n\right]  ^{2}$, so that we can write $c$ in
the form $c=\left(  i,j\right)  $ for some $i\in\left[  n\right]  $ and
$j\in\left[  n\right]  $. Consider these $i$ and $j$.
\par
We have $\left(  i,j\right)  =c\in\left(  s_{k}\times s_{k}\right)
^{-1}\left(  \operatorname*{Inv}\left(  \mathbf{a}\right)  \setminus\left\{
\left(  k,k+1\right)  \right\}  \right)  $, so that $\left(  s_{k}\times
s_{k}\right)  \left(  \left(  i,j\right)  \right)  \in\operatorname*{Inv}%
\left(  \mathbf{a}\right)  \setminus\left\{  \left(  k,k+1\right)  \right\}
$. Since $\left(  s_{k}\times s_{k}\right)  \left(  \left(  i,j\right)
\right)  =\left(  s_{k}\left(  i\right)  ,s_{k}\left(  j\right)  \right)  $
(by the definition of $\left(  s_{k}\times s_{k}\right)  $), this rewrites as
$\left(  s_{k}\left(  i\right)  ,s_{k}\left(  j\right)  \right)
\in\operatorname*{Inv}\left(  \mathbf{a}\right)  \setminus\left\{  \left(
k,k+1\right)  \right\}  $. In other words, $\left(  s_{k}\left(  i\right)
,s_{k}\left(  j\right)  \right)  \in\operatorname*{Inv}\left(  \mathbf{a}%
\right)  $ and $\left(  s_{k}\left(  i\right)  ,s_{k}\left(  j\right)
\right)  \neq\left(  k,k+1\right)  $.
\par
We have $\left(  s_{k}\left(  i\right)  ,s_{k}\left(  j\right)  \right)
\in\operatorname*{Inv}\left(  \mathbf{a}\right)  =\left\{  \left(  u,v\right)
\in\left[  n\right]  ^{2}\ \mid\ u<v\text{ and }a_{u}>a_{v}\right\}  $ (by
(\ref{pf.prop.sorting.a.Inv.2})). In other words, $\left(  s_{k}\left(
i\right)  ,s_{k}\left(  j\right)  \right)  $ is an element of $\left[
n\right]  ^{2}$ and satisfies $s_{k}\left(  i\right)  <s_{k}\left(  j\right)
$ and $a_{s_{k}\left(  i\right)  }>a_{s_{k}\left(  j\right)  }$. We have
$s_{k}\left(  i\right)  <s_{k}\left(  j\right)  $ and $\left(  s_{k}\left(
i\right)  ,s_{k}\left(  j\right)  \right)  \neq\left(  k,k+1\right)  $. Thus,
we can apply (\ref{pf.prop.sorting.a.2.5}) to $u=s_{k}\left(  i\right)  $ and
$v=s_{k}\left(  j\right)  $. We thus conclude $s_{k}\left(  s_{k}\left(
i\right)  \right)  <s_{k}\left(  s_{k}\left(  j\right)  \right)  $. In other
words, $i<j$ (since $s_{k}\left(  s_{k}\left(  i\right)  \right)
=\underbrace{s_{k}^{2}}_{=\operatorname*{id}}\left(  i\right)
=\operatorname*{id}\left(  i\right)  =i$ and $s_{k}\left(  s_{k}\left(
j\right)  \right)  =\underbrace{s_{k}^{2}}_{=\operatorname*{id}}\left(
j\right)  =\operatorname*{id}\left(  j\right)  =j$).
\par
Now, we know that $\left(  i,j\right)  $ is an element $\left(  u,v\right)  $
of $\left[  n\right]  ^{2}$ satisfying $u<v$ and $a_{s_{k}\left(  u\right)
}>a_{s_{k}\left(  v\right)  }$ (since $i<j$ and $a_{s_{k}\left(  i\right)
}>a_{s_{k}\left(  j\right)  }$). In other words,%
\[
\left(  i,j\right)  \in\left\{  \left(  u,v\right)  \in\left[  n\right]
^{2}\ \mid\ u<v\text{ and }a_{s_{k}\left(  u\right)  }>a_{s_{k}\left(
v\right)  }\right\}  =\operatorname*{Inv}\left(  \mathbf{a}\circ s_{k}\right)
\]
(by (\ref{pf.prop.sorting.a.fact1})).
\par
Next, let us assume (for the sake of contradiction) that $\left(  i,j\right)
=\left(  k,k+1\right)  $. Thus, $i=k$ and $j=k+1$. Thus, $s_{k}\left(
\underbrace{i}_{=k}\right)  =s_{k}\left(  k\right)  =k+1>k=s_{k}\left(
\underbrace{k+1}_{=j}\right)  =s_{k}\left(  j\right)  $; but this contradicts
$s_{k}\left(  i\right)  <s_{k}\left(  j\right)  $. This contradiction shows
that our assumption (that $\left(  i,j\right)  =\left(  k,k+1\right)  $) was
wrong. Hence, we cannot have $\left(  i,j\right)  =\left(  k,k+1\right)  $. In
other words, we must have $\left(  i,j\right)  \neq\left(  k,k+1\right)  $.
Combined with $\left(  i,j\right)  \in\operatorname*{Inv}\left(
\mathbf{a}\circ s_{k}\right)  $, this yields $\left(  i,j\right)
\in\operatorname*{Inv}\left(  \mathbf{a}\circ s_{k}\right)  \setminus\left\{
\left(  k,k+1\right)  \right\}  $. Thus, $c=\left(  i,j\right)  \in
\operatorname*{Inv}\left(  \mathbf{a}\circ s_{k}\right)  \setminus\left\{
\left(  k,k+1\right)  \right\}  $.
\par
Now, let us forget that we fixed $c$. We thus have shown that $c\in
\operatorname*{Inv}\left(  \mathbf{a}\circ s_{k}\right)  \setminus\left\{
\left(  k,k+1\right)  \right\}  $ for every $c\in\left(  s_{k}\times
s_{k}\right)  ^{-1}\left(  \operatorname*{Inv}\left(  \mathbf{a}\right)
\setminus\left\{  \left(  k,k+1\right)  \right\}  \right)  $. In other words,%
\[
\left(  s_{k}\times s_{k}\right)  ^{-1}\left(  \operatorname*{Inv}\left(
\mathbf{a}\right)  \setminus\left\{  \left(  k,k+1\right)  \right\}  \right)
\subseteq\operatorname*{Inv}\left(  \mathbf{a}\circ s_{k}\right)
\setminus\left\{  \left(  k,k+1\right)  \right\}  .
\]
This proves (\ref{pf.prop.sorting.a.fact2}).}.

\item If $\mathbf{a}$ is an $n$-tuple of integers, and if $k\in\left\{
1,2,\ldots,n-1\right\}  $, then%
\begin{equation}
\operatorname*{Inv}\left(  \mathbf{a}\circ s_{k}\right)  \setminus\left\{
\left(  k,k+1\right)  \right\}  \subseteq\left(  s_{k}\times s_{k}\right)
^{-1}\left(  \operatorname*{Inv}\left(  \mathbf{a}\right)  \setminus\left\{
\left(  k,k+1\right)  \right\}  \right)  \label{pf.prop.sorting.a.fact3}%
\end{equation}
\footnote{\textit{Proof of (\ref{pf.prop.sorting.a.fact3}):} Let $\mathbf{a}$
be an $n$-tuple of integers, and let $k\in\left\{  1,2,\ldots,n-1\right\}  $.
Write the $n$-tuple $\mathbf{a}$ in the form $\mathbf{a}=\left(  a_{1}%
,a_{2},\ldots,a_{n}\right)  $ for some integers $a_{1},a_{2},\ldots,a_{n}$.
\par
Let $c\in\operatorname*{Inv}\left(  \mathbf{a}\circ s_{k}\right)
\setminus\left\{  \left(  k,k+1\right)  \right\}  $. Thus, $c\in
\operatorname*{Inv}\left(  \mathbf{a}\circ s_{k}\right)  \setminus\left\{
\left(  k,k+1\right)  \right\}  \subseteq\operatorname*{Inv}\left(
\mathbf{a}\circ s_{k}\right)  \subseteq\left[  n\right]  ^{2}$, so that we can
write $c$ in the form $c=\left(  i,j\right)  $ for some $i\in\left[  n\right]
$ and $j\in\left[  n\right]  $. Consider these $i$ and $j$.
\par
We have $\left(  i,j\right)  =c\in\operatorname*{Inv}\left(  \mathbf{a}\circ
s_{k}\right)  \setminus\left\{  \left(  k,k+1\right)  \right\}  $. In other
words, $\left(  i,j\right)  \in\operatorname*{Inv}\left(  \mathbf{a}\circ
s_{k}\right)  $ and $\left(  i,j\right)  \neq\left(  k,k+1\right)  $.
\par
We have $\left(  i,j\right)  \in\operatorname*{Inv}\left(  \mathbf{a}\circ
s_{k}\right)  =\left\{  \left(  u,v\right)  \in\left[  n\right]  ^{2}%
\ \mid\ u<v\text{ and }a_{s_{k}\left(  u\right)  }>a_{s_{k}\left(  v\right)
}\right\}  $ (by (\ref{pf.prop.sorting.a.fact1})). In other words, $\left(
i,j\right)  $ is an element of $\left[  n\right]  ^{2}$ and satisfies $i<j$
and $a_{s_{k}\left(  i\right)  }>a_{s_{k}\left(  j\right)  }$. Applying
(\ref{pf.prop.sorting.a.2.5}) to $u=i$ and $v=j$, we obtain $s_{k}\left(
i\right)  <s_{k}\left(  j\right)  $ (since $i<j$ and $\left(  i,j\right)
\neq\left(  k,k+1\right)  $).
\par
The pair $\left(  s_{k}\left(  i\right)  ,s_{k}\left(  j\right)  \right)  $ is
a pair $\left(  u,v\right)  \in\left[  n\right]  ^{2}$ satisfying $u<v$ and
$a_{u}>a_{v}$ (since $s_{k}\left(  i\right)  <s_{k}\left(  j\right)  $ and
$a_{s_{k}\left(  i\right)  }>a_{s_{k}\left(  j\right)  }$). In other words,
$\left(  s_{k}\left(  i\right)  ,s_{k}\left(  j\right)  \right)  \in\left\{
\left(  u,v\right)  \in\left[  n\right]  ^{2}\ \mid\ u<v\text{ and }%
a_{u}>a_{v}\right\}  =\operatorname*{Inv}\left(  \mathbf{a}\right)  $ (by
(\ref{pf.prop.sorting.a.Inv.2})).
\par
Next, let us assume (for the sake of contradiction) that $\left(  s_{k}\left(
i\right)  ,s_{k}\left(  j\right)  \right)  =\left(  k,k+1\right)  $. Thus,
$s_{k}\left(  i\right)  =k$ and $s_{k}\left(  j\right)  =k+1$. Hence,
$k+1=s_{k}\left(  \underbrace{k}_{=s_{k}\left(  i\right)  }\right)
=s_{k}\left(  s_{k}\left(  i\right)  \right)  =\underbrace{s_{k}^{2}%
}_{=\operatorname*{id}}\left(  i\right)  =i$ and $k=s_{k}\left(
\underbrace{k+1}_{=s_{k}\left(  j\right)  }\right)  =s_{k}\left(  s_{k}\left(
j\right)  \right)  =\underbrace{s_{k}^{2}}_{=\operatorname*{id}}\left(
j\right)  =\operatorname*{id}\left(  j\right)  =j$, so that $i=k+1>k=j$. This
contradicts $i<j$. This contradiction shows that our assumption (that $\left(
s_{k}\left(  i\right)  ,s_{k}\left(  j\right)  \right)  =\left(  k,k+1\right)
$) was wrong. Hence, we cannot have $\left(  s_{k}\left(  i\right)
,s_{k}\left(  j\right)  \right)  =\left(  k,k+1\right)  $. In other words, we
must have $\left(  s_{k}\left(  i\right)  ,s_{k}\left(  j\right)  \right)
\neq\left(  k,k+1\right)  $. Combined with $\left(  s_{k}\left(  i\right)
,s_{k}\left(  j\right)  \right)  \in\operatorname*{Inv}\left(  \mathbf{a}%
\right)  $, this yields $\left(  s_{k}\left(  i\right)  ,s_{k}\left(
j\right)  \right)  \in\operatorname*{Inv}\left(  \mathbf{a}\right)
\setminus\left\{  \left(  k,k+1\right)  \right\}  $.
\par
The definition of $s_{k}\times s_{k}$ yields $\left(  s_{k}\times
s_{k}\right)  \left(  \left(  i,j\right)  \right)  =\left(  s_{k}\left(
i\right)  ,s_{k}\left(  j\right)  \right)  \in\operatorname*{Inv}\left(
\mathbf{a}\right)  \setminus\left\{  \left(  k,k+1\right)  \right\}  $. In
other words, $\left(  i,j\right)  \in\left(  s_{k}\times s_{k}\right)
^{-1}\left(  \operatorname*{Inv}\left(  \mathbf{a}\right)  \setminus\left\{
\left(  k,k+1\right)  \right\}  \right)  $. Hence,%
\[
c=\left(  i,j\right)  \in\left(  s_{k}\times s_{k}\right)  ^{-1}\left(
\operatorname*{Inv}\left(  \mathbf{a}\right)  \setminus\left\{  \left(
k,k+1\right)  \right\}  \right)  .
\]
\par
Now, let us forget that we fixed $c$. We thus have shown that $c\in\left(
s_{k}\times s_{k}\right)  ^{-1}\left(  \operatorname*{Inv}\left(
\mathbf{a}\right)  \setminus\left\{  \left(  k,k+1\right)  \right\}  \right)
$ for every $c\in\operatorname*{Inv}\left(  \mathbf{a}\circ s_{k}\right)
\setminus\left\{  \left(  k,k+1\right)  \right\}  $. In other words,%
\[
\operatorname*{Inv}\left(  \mathbf{a}\circ s_{k}\right)  \setminus\left\{
\left(  k,k+1\right)  \right\}  \subseteq\left(  s_{k}\times s_{k}\right)
^{-1}\left(  \operatorname*{Inv}\left(  \mathbf{a}\right)  \setminus\left\{
\left(  k,k+1\right)  \right\}  \right)  .
\]
This proves (\ref{pf.prop.sorting.a.fact3}).}.

\item If $\mathbf{a}$ is an $n$-tuple of integers, and if $k\in\left\{
1,2,\ldots,n-1\right\}  $, then%
\begin{equation}
\left(  s_{k}\times s_{k}\right)  ^{-1}\left(  \operatorname*{Inv}\left(
\mathbf{a}\right)  \setminus\left\{  \left(  k,k+1\right)  \right\}  \right)
=\operatorname*{Inv}\left(  \mathbf{a}\circ s_{k}\right)  \setminus\left\{
\left(  k,k+1\right)  \right\}  \label{pf.prop.sorting.a.fact4}%
\end{equation}
\footnote{\textit{Proof of (\ref{pf.prop.sorting.a.fact4}):} Let $\mathbf{a}$
be an $n$-tuple of integers, and let $k\in\left\{  1,2,\ldots,n-1\right\}  $.
Then, combining (\ref{pf.prop.sorting.a.fact2}) with
(\ref{pf.prop.sorting.a.fact3}), we obtain $\left(  s_{k}\times s_{k}\right)
^{-1}\left(  \operatorname*{Inv}\left(  \mathbf{a}\right)  \setminus\left\{
\left(  k,k+1\right)  \right\}  \right)  =\operatorname*{Inv}\left(
\mathbf{a}\circ s_{k}\right)  \setminus\left\{  \left(  k,k+1\right)
\right\}  $. This proves (\ref{pf.prop.sorting.a.fact4}).}.

\item If $\mathbf{a}=\left(  a_{1},a_{2},\ldots,a_{n}\right)  $ is an
$n$-tuple of integers, and if $k\in\left\{  1,2,\ldots,n-1\right\}  $ is such
that $a_{k}>a_{k+1}$, then%
\begin{equation}
\ell\left(  \mathbf{a}\circ s_{k}\right)  =\ell\left(  \mathbf{a}\right)  -1
\label{pf.prop.sorting.a.fact5}%
\end{equation}
\footnote{\textit{Proof of (\ref{pf.prop.sorting.a.fact5}):} Let
$\mathbf{a}=\left(  a_{1},a_{2},\ldots,a_{n}\right)  $ be an $n$-tuple of
integers, and let $k\in\left\{  1,2,\ldots,n-1\right\}  $ be such that
$a_{k}>a_{k+1}$.
\par
The definition of $\ell\left(  \mathbf{a}\right)  $ yields $\ell\left(
\mathbf{a}\right)  =\left\vert \operatorname*{Inv}\left(  \mathbf{a}\right)
\right\vert $.
\par
The pair $\left(  k,k+1\right)  $ is a pair $\left(  u,v\right)  \in\left[
n\right]  ^{2}$ such that $u<v$ and $a_{u}>a_{v}$ (since $k<k+1$ and
$a_{k}>a_{k+1}$). In other words,%
\[
\left(  k,k+1\right)  \in\left\{  \left(  u,v\right)  \in\left[  n\right]
^{2}\ \mid\ u<v\text{ and }a_{u}>a_{v}\right\}  =\operatorname*{Inv}\left(
\mathbf{a}\right)  \ \ \ \ \ \ \ \ \ \ \left(  \text{by
(\ref{pf.prop.sorting.a.Inv.2})}\right)  .
\]
Hence, $\left\vert \operatorname*{Inv}\left(  \mathbf{a}\right)
\setminus\left\{  \left(  k,k+1\right)  \right\}  \right\vert
=\underbrace{\left\vert \operatorname*{Inv}\left(  \mathbf{a}\right)
\right\vert }_{=\ell\left(  \mathbf{a}\right)  }-1=\ell\left(  \mathbf{a}%
\right)  -1$. The map $s_{k}\times s_{k}$ is a bijection, and thus we have
$\left\vert \left(  s_{k}\times s_{k}\right)  ^{-1}\left(  X\right)
\right\vert =\left\vert X\right\vert $ for every subset $X$ of $\left[
n\right]  ^{2}$. Applying this to $X=\operatorname*{Inv}\left(  \mathbf{a}%
\right)  \setminus\left\{  \left(  k,k+1\right)  \right\}  $, we obtain%
\begin{equation}
\left\vert \left(  s_{k}\times s_{k}\right)  ^{-1}\left(  \operatorname*{Inv}%
\left(  \mathbf{a}\right)  \setminus\left\{  \left(  k,k+1\right)  \right\}
\right)  \right\vert =\left\vert \operatorname*{Inv}\left(  \mathbf{a}\right)
\setminus\left\{  \left(  k,k+1\right)  \right\}  \right\vert =\ell\left(
\mathbf{a}\right)  -1. \label{pf.prop.sorting.a.fact5.pf.1}%
\end{equation}
\par
On the other hand, let us assume (for the sake of contradiction) that $\left(
k,k+1\right)  \in\operatorname*{Inv}\left(  \mathbf{a}\circ s_{k}\right)  $.
Thus,%
\[
\left(  k,k+1\right)  \in\operatorname*{Inv}\left(  \mathbf{a}\circ
s_{k}\right)  =\left\{  \left(  u,v\right)  \in\left[  n\right]  ^{2}%
\ \mid\ u<v\text{ and }a_{s_{k}\left(  u\right)  }>a_{s_{k}\left(  v\right)
}\right\}
\]
(by (\ref{pf.prop.sorting.a.fact1})). In other words, $\left(  k,k+1\right)  $
is a pair $\left(  u,v\right)  \in\left[  n\right]  ^{2}$ such that $u<v$ and
$a_{s_{k}\left(  u\right)  }>a_{s_{k}\left(  v\right)  }$. In other words,
$k<k+1$ and $a_{s_{k}\left(  k\right)  }>a_{s_{k}\left(  k+1\right)  }$. Now,
$a_{s_{k}\left(  k\right)  }>a_{s_{k}\left(  k+1\right)  }$. In other words,
$a_{k+1}>a_{k}$ (since $s_{k}\left(  k\right)  =k+1$ and $s_{k}\left(
k+1\right)  =k$). This contradicts $a_{k}>a_{k+1}$. This contradiction shows
that our assumption (that $\left(  k,k+1\right)  \in\operatorname*{Inv}\left(
\mathbf{a}\circ s_{k}\right)  $) was wrong. Hence, we cannot have $\left(
k,k+1\right)  \in\operatorname*{Inv}\left(  \mathbf{a}\circ s_{k}\right)  $.
We thus have $\left(  k,k+1\right)  \notin\operatorname*{Inv}\left(
\mathbf{a}\circ s_{k}\right)  $. Hence, $\operatorname*{Inv}\left(
\mathbf{a}\circ s_{k}\right)  \setminus\left\{  \left(  k,k+1\right)
\right\}  =\operatorname*{Inv}\left(  \mathbf{a}\circ s_{k}\right)  $. Now,
(\ref{pf.prop.sorting.a.fact4}) becomes%
\[
\left(  s_{k}\times s_{k}\right)  ^{-1}\left(  \operatorname*{Inv}\left(
\mathbf{a}\right)  \setminus\left\{  \left(  k,k+1\right)  \right\}  \right)
=\operatorname*{Inv}\left(  \mathbf{a}\circ s_{k}\right)  \setminus\left\{
\left(  k,k+1\right)  \right\}  =\operatorname*{Inv}\left(  \mathbf{a}\circ
s_{k}\right)  .
\]
Therefore, (\ref{pf.prop.sorting.a.fact5.pf.1}) rewrites as $\left\vert
\operatorname*{Inv}\left(  \mathbf{a}\circ s_{k}\right)  \right\vert
=\ell\left(  \mathbf{a}\right)  -1$.
\par
But the definition of $\ell\left(  \mathbf{a}\circ s_{k}\right)  $ yields
$\ell\left(  \mathbf{a}\circ s_{k}\right)  =\left\vert \operatorname*{Inv}%
\left(  \mathbf{a}\circ s_{k}\right)  \right\vert =\ell\left(  \mathbf{a}%
\right)  -1$. This proves (\ref{pf.prop.sorting.a.fact5}).}.

\item If $\mathbf{a}=\left(  a_{1},a_{2},\ldots,a_{n}\right)  $ is an
$n$-tuple of integers satisfying $\ell\left(  \mathbf{a}\right)  =0$, then%
\begin{equation}
a_{1}\leq a_{2}\leq\cdots\leq a_{n} \label{pf.prop.sorting.a.fact7}%
\end{equation}
\footnote{\textit{Proof of (\ref{pf.prop.sorting.a.fact7}):} Let
$\mathbf{a}=\left(  a_{1},a_{2},\ldots,a_{n}\right)  $ be an $n$-tuple of
integers such that $\ell\left(  \mathbf{a}\right)  =0$.
\par
We assume (for the sake of contradiction) that there exists some $k\in\left\{
1,2,\ldots,n-1\right\}  $ such that $a_{k}>a_{k+1}$. Consider this $k$. Then,
$\left(  k,k+1\right)  $ is an element $\left(  u,v\right)  \in\left[
n\right]  ^{2}$ satisfying $u<v$ and $a_{u}>a_{v}$ (since $k<k+1$ and
$a_{k}>a_{k+1}$). In other words,%
\[
\left(  k,k+1\right)  \in\left\{  \left(  u,v\right)  \in\left[  n\right]
^{2}\ \mid\ u<v\text{ and }a_{u}>a_{v}\right\}  =\operatorname*{Inv}\left(
\mathbf{a}\right)
\]
(by (\ref{pf.prop.sorting.a.Inv.2})). Hence, the set $\operatorname*{Inv}%
\left(  \mathbf{a}\right)  $ contains at least one element (namely, $\left(
k,k+1\right)  $). In other words, $\left\vert \operatorname*{Inv}\left(
\mathbf{a}\right)  \right\vert \geq1$.
\par
But the definition of $\ell\left(  \mathbf{a}\right)  $ yields $\ell\left(
\mathbf{a}\right)  =\left\vert \operatorname*{Inv}\left(  \mathbf{a}\right)
\right\vert \geq1$. This contradicts $\ell\left(  \mathbf{a}\right)  =0$. This
contradiction shows that our assumption (that there exists some $k\in\left\{
1,2,\ldots,n-1\right\}  $ such that $a_{k}>a_{k+1}$) was wrong.
\par
Therefore, there exists no $k\in\left\{  1,2,\ldots,n-1\right\}  $ such that
$a_{k}>a_{k+1}$. In other words, every $k\in\left\{  1,2,\ldots,n-1\right\}  $
satisfies $a_{k}\leq a_{k+1}$. In other words, $a_{1}\leq a_{2}\leq\cdots\leq
a_{n}$. This proves (\ref{pf.prop.sorting.a.fact7}).}.

\item If $\mathbf{a}=\left(  a_{1},a_{2},\ldots,a_{n}\right)  $ is an
$n$-tuple of integers satisfying $a_{1}\leq a_{2}\leq\cdots\leq a_{n}$, then%
\begin{equation}
\ell\left(  \mathbf{a}\right)  =0 \label{pf.prop.sorting.a.fact8}%
\end{equation}
\footnote{\textit{Proof of (\ref{pf.prop.sorting.a.fact8}):} Let
$\mathbf{a}=\left(  a_{1},a_{2},\ldots,a_{n}\right)  $ be an $n$-tuple of
integers satisfying $a_{1}\leq a_{2}\leq\cdots\leq a_{n}$. We must prove that
$\ell\left(  \mathbf{a}\right)  =0$.
\par
Indeed, assume the contrary. Thus, $\ell\left(  \mathbf{a}\right)  \neq0$. But
the definition of $\ell\left(  \mathbf{a}\right)  $ yields $\ell\left(
\mathbf{a}\right)  =\left\vert \operatorname*{Inv}\left(  \mathbf{a}\right)
\right\vert $, so that $\left\vert \operatorname*{Inv}\left(  \mathbf{a}%
\right)  \right\vert =\ell\left(  \mathbf{a}\right)  \neq0$. Hence, the set
$\operatorname*{Inv}\left(  \mathbf{a}\right)  $ is nonempty. In other words,
there exists some $c\in\operatorname*{Inv}\left(  \mathbf{a}\right)  $.
Consider this $c$.
\par
Recall that $a_{1}\leq a_{2}\leq\cdots\leq a_{n}$. In other words,%
\begin{equation}
a_{u}\leq a_{v} \label{pf.prop.sorting.a.fact8.pf.leq}%
\end{equation}
for any $u\in\left[  n\right]  $ and $v\in\left[  n\right]  $ satisfying
$u<v$.
\par
But%
\[
c\in\operatorname*{Inv}\left(  \mathbf{a}\right)  =\left\{  \left(
u,v\right)  \in\left[  n\right]  ^{2}\ \mid\ u<v\text{ and }a_{u}%
>a_{v}\right\}
\]
(by (\ref{pf.prop.sorting.a.Inv.2})). In other words, $c$ has the form
$c=\left(  u,v\right)  $ for some $\left(  u,v\right)  \in\left[  n\right]
^{2}$ satisfying $u<v$ and $a_{u}>a_{v}$. Consider this $\left(  u,v\right)
\in\left[  n\right]  ^{2}$. We have $u<v$, and thus $a_{u}\leq a_{v}$ (by
(\ref{pf.prop.sorting.a.fact8.pf.leq})). This contradicts $a_{u}>a_{v}$. This
contradiction proves that our assumption was wrong.
\par
Hence, we have shown that $\ell\left(  \mathbf{a}\right)  =0$. This proves
(\ref{pf.prop.sorting.a.fact8}).}.

\item If $\mathbf{a}$ is an $n$-tuple of integers, then%
\begin{equation}
\text{there exists a }\sigma\in S_{n}\text{ such that }\ell\left(
\mathbf{a}\circ\sigma\right)  =0 \label{pf.prop.sorting.a.fact6}%
\end{equation}
\footnote{\textit{Proof of (\ref{pf.prop.sorting.a.fact6}):} We shall prove
(\ref{pf.prop.sorting.a.fact6}) by induction over $\ell\left(  \mathbf{a}%
\right)  $:
\par
\textit{Induction base:} If $\mathbf{a}$ is an $n$-tuple of integers
satisfying $\ell\left(  \mathbf{a}\right)  =0$, then we have $\ell\left(
\underbrace{\mathbf{a}\circ\operatorname*{id}}_{\substack{=\mathbf{a}%
\\\text{(by (\ref{pf.prop.sorting.a.id}))}}}\right)  =\ell\left(
\mathbf{a}\right)  =0$. Hence, if $\mathbf{a}$ is an $n$-tuple of integers
satisfying $\ell\left(  \mathbf{a}\right)  =0$, then there exists a $\sigma\in
S_{n}$ such that $\ell\left(  \mathbf{a}\circ\sigma\right)  =0$ (namely,
$\sigma=\operatorname*{id}$). In other words, (\ref{pf.prop.sorting.a.fact6})
holds for $\ell\left(  \mathbf{a}\right)  =0$. This completes the induction
base.
\par
\textit{Induction step:} Let $L$ be a positive integer. Assume that
(\ref{pf.prop.sorting.a.fact6}) holds for $\ell\left(  \mathbf{a}\right)
=L-1$. We need to prove that (\ref{pf.prop.sorting.a.fact6}) holds for
$\ell\left(  \mathbf{a}\right)  =L$.
\par
We have assumed that (\ref{pf.prop.sorting.a.fact6}) holds for $\ell\left(
\mathbf{a}\right)  =L-1$. In other words, if $\mathbf{a}$ is an $n$-tuple of
integers satisfying $\ell\left(  \mathbf{a}\right)  =L-1$, then%
\begin{equation}
\text{there exists a }\sigma\in S_{n}\text{ such that }\ell\left(
\mathbf{a}\circ\sigma\right)  =0. \label{pf.prop.sorting.a.fact6.pf.hyp}%
\end{equation}
\par
Now, let $\mathbf{a}$ be an $n$-tuple of integers satisfying $\ell\left(
\mathbf{a}\right)  =L$. We shall show that there exists a $\sigma\in S_{n}$
such that $\ell\left(  \mathbf{a}\circ\sigma\right)  =0$.
\par
Let us first prove that there exists some $k\in\left\{  1,2,\ldots
,n-1\right\}  $ such that $a_{k}>a_{k+1}$. In fact, assume the contrary. Thus,
there exists no $k\in\left\{  1,2,\ldots,n-1\right\}  $ such that
$a_{k}>a_{k+1}$. In other words, every $k\in\left\{  1,2,\ldots,n-1\right\}  $
satisfies $a_{k}\leq a_{k+1}$. In other words, $a_{1}\leq a_{2}\leq\cdots\leq
a_{n}$. Thus, (\ref{pf.prop.sorting.a.fact8}) shows that $\ell\left(
\mathbf{a}\right)  =0$. Hence, $0=\ell\left(  \mathbf{a}\right)  =L>0$ (since
$L$ is positive). This is absurd. This contradiction proves that our
assumption was wrong.
\par
Hence, we have proven that there exists some $k\in\left\{  1,2,\ldots
,n-1\right\}  $ such that $a_{k}>a_{k+1}$. Consider this $k$. Then,
(\ref{pf.prop.sorting.a.fact5}) shows that $\ell\left(  \mathbf{a}\circ
s_{k}\right)  =\underbrace{\ell\left(  \mathbf{a}\right)  }_{=L}-1=L-1$.
Therefore, (\ref{pf.prop.sorting.a.fact6.pf.hyp}) (applied to $\mathbf{a}\circ
s_{k}$ instead of $\mathbf{a}$) shows that there exists a $\sigma\in S_{n}$
such that $\ell\left(  \left(  \mathbf{a}\circ s_{k}\right)  \circ
\sigma\right)  =0$. Let $\tau$ be such a $\sigma$. Thus, $\tau$ is an element
of $S_{n}$ and satisfies $\ell\left(  \left(  \mathbf{a}\circ s_{k}\right)
\circ\tau\right)  =0$.
\par
Now, (\ref{pf.prop.sorting.a.sigmatau}) (applied to $s_{k}$ instead of
$\sigma$) yields $\mathbf{a}\circ\left(  s_{k}\circ\tau\right)  =\left(
\mathbf{a}\circ s_{k}\right)  \circ\tau$. Hence, $\ell\left(  \mathbf{a}%
\circ\left(  s_{k}\circ\tau\right)  \right)  =\ell\left(  \left(
\mathbf{a}\circ s_{k}\right)  \circ\tau\right)  =0$. Thus, there exists a
$\sigma\in S_{n}$ such that $\ell\left(  \mathbf{a}\circ\sigma\right)  =0$
(namely, $\sigma=s_{k}\circ\tau$).
\par
Now, let us forget that we fixed $\mathbf{a}$. We thus have shown that if
$\mathbf{a}$ is an $n$-tuple of integers satisfying $\ell\left(
\mathbf{a}\right)  =L$, then there exists a $\sigma\in S_{n}$ such that
$\ell\left(  \mathbf{a}\circ\sigma\right)  =0$. In other words,
(\ref{pf.prop.sorting.a.fact6}) holds for $\ell\left(  \mathbf{a}\right)  =L$.
This completes the induction step. Thus, the induction proof of
(\ref{pf.prop.sorting.a.fact6}) is complete.}.
\end{itemize}

Now, let us prove Proposition \ref{prop.sorting} \textbf{(a)}. Let
$a_{1},a_{2},\ldots,a_{n}$ be $n$ integers. Let $\mathbf{a}$ be the $n$-tuple
$\left(  a_{1},a_{2},\ldots,a_{n}\right)  $. Then, there exists a $\sigma\in
S_{n}$ such that $\ell\left(  \mathbf{a}\circ\sigma\right)  =0$ (by
(\ref{pf.prop.sorting.a.fact6})). Consider this $\sigma$. The definition of
$\mathbf{a}\circ\sigma$ shows that $\mathbf{a}\circ\sigma=\left(
a_{\sigma\left(  1\right)  },a_{\sigma\left(  2\right)  },\ldots
,a_{\sigma\left(  n\right)  }\right)  $ (since $\mathbf{a}=\left(  a_{1}%
,a_{2},\ldots,a_{n}\right)  $). Therefore, (\ref{pf.prop.sorting.a.fact7})
(applied to $\mathbf{a}\circ\sigma$ and $a_{\sigma\left(  i\right)  }$ instead
of $\mathbf{a}$ and $a_{i}$) yields $a_{\sigma\left(  1\right)  }\leq
a_{\sigma\left(  2\right)  }\leq\cdots\leq a_{\sigma\left(  n\right)  }$.

Let us now forget that we defined $\sigma$. We thus have constructed a
$\sigma\in S_{n}$ satisfying $a_{\sigma\left(  1\right)  }\leq a_{\sigma
\left(  2\right)  }\leq\cdots\leq a_{\sigma\left(  n\right)  }$. Therefore,
there exists a permutation $\sigma\in S_{n}$ such that $a_{\sigma\left(
1\right)  }\leq a_{\sigma\left(  2\right)  }\leq\cdots\leq a_{\sigma\left(
n\right)  }$. This proves Proposition \ref{prop.sorting} \textbf{(a)}.

\textbf{(b)} Let $\sigma\in S_{n}$ be such that $a_{\sigma\left(  1\right)
}\leq a_{\sigma\left(  2\right)  }\leq\cdots\leq a_{\sigma\left(  n\right)  }%
$. Let $i\in\left\{  1,2,\ldots,n\right\}  $. We shall now show that%
\begin{align}
&  a_{\sigma\left(  i\right)  }\nonumber\\
&  =\min\left\{  x\in\mathbb{Z}\ \mid\ \text{at least }i\text{ elements }%
j\in\left\{  1,2,\ldots,n\right\}  \text{ satisfy }a_{j}\leq x\right\}
.\ \ \ \ \ \ \label{pf.prop.sorting.bclaim}%
\end{align}
(This statement includes the tacit claim that the right hand side of
(\ref{pf.prop.sorting.bclaim}) is well-defined.)

[\textit{Proof of (\ref{pf.prop.sorting.bclaim}):} Define a subset $X$ of
$\mathbb{Z}$ by%
\[
X=\left\{  x\in\mathbb{Z}\ \mid\ \text{at least }i\text{ elements }%
j\in\left\{  1,2,\ldots,n\right\}  \text{ satisfy }a_{j}\leq x\right\}  .
\]
We shall show that $a_{\sigma\left(  i\right)  }=\min X$.

We have $a_{\sigma\left(  1\right)  }\leq a_{\sigma\left(  2\right)  }%
\leq\cdots\leq a_{\sigma\left(  n\right)  }$. In other words,%
\begin{equation}
a_{\sigma\left(  u\right)  }\leq a_{\sigma\left(  v\right)  }
\label{pf.prop.sorting.bclaim.pf.1}%
\end{equation}
for any two elements $u$ and $v$ of $\left\{  1,2,\ldots,n\right\}  $
satisfying $u\leq v$.

We have $\sigma\in S_{n}$. Thus, $\sigma$ is a permutation, thus a bijective
map, thus an injective map.

Every $k\in\left\{  1,2,\ldots,i\right\}  $ satisfies $\sigma\left(  k\right)
\in\left\{  j\in\left\{  1,2,\ldots,n\right\}  \ \mid\ a_{j}\leq
a_{\sigma\left(  i\right)  }\right\}  $\ \ \ \ \footnote{\textit{Proof.} Let
$k\in\left\{  1,2,\ldots,i\right\}  $. Then, $k\leq i$ and thus $a_{\sigma
\left(  k\right)  }\leq a_{\sigma\left(  i\right)  }$ (by
(\ref{pf.prop.sorting.bclaim.pf.1}), applied to $u=k$ and $v=i$). Hence,
$\sigma\left(  k\right)  $ is an element $j$ of $\left\{  1,2,\ldots
,n\right\}  $ satisfying $a_{j}\leq a_{\sigma\left(  i\right)  }$ (because
$a_{\sigma\left(  k\right)  }\leq a_{\sigma\left(  i\right)  }$). In other
words, $\sigma\left(  k\right)  \in\left\{  j\in\left\{  1,2,\ldots,n\right\}
\ \mid\ a_{j}\leq a_{\sigma\left(  i\right)  }\right\}  $, qed.}. In other
words,%
\[
\left\{  \sigma\left(  k\right)  \ \mid\ k\in\left\{  1,2,\ldots,i\right\}
\right\}  \subseteq\left\{  j\in\left\{  1,2,\ldots,n\right\}  \ \mid
\ a_{j}\leq a_{\sigma\left(  i\right)  }\right\}  .
\]
Hence,%
\[
\left\vert \left\{  \sigma\left(  k\right)  \ \mid\ k\in\left\{
1,2,\ldots,i\right\}  \right\}  \right\vert \leq\left\vert \left\{
j\in\left\{  1,2,\ldots,n\right\}  \ \mid\ a_{j}\leq a_{\sigma\left(
i\right)  }\right\}  \right\vert ,
\]
so that
\begin{align*}
\left\vert \left\{  j\in\left\{  1,2,\ldots,n\right\}  \ \mid\ a_{j}\leq
a_{\sigma\left(  i\right)  }\right\}  \right\vert  &  \geq\left\vert
\underbrace{\left\{  \sigma\left(  k\right)  \ \mid\ k\in\left\{
1,2,\ldots,i\right\}  \right\}  }_{=\sigma\left(  \left\{  1,2,\ldots
,i\right\}  \right)  }\right\vert \\
&  =\left\vert \sigma\left(  \left\{  1,2,\ldots,i\right\}  \right)
\right\vert =\left\vert \left\{  1,2,\ldots,i\right\}  \right\vert \\
&  \ \ \ \ \ \ \ \ \ \ \left(  \text{since the map }\sigma\text{ is
injective}\right) \\
&  =i.
\end{align*}
In other words, at least $i$ elements $j\in\left\{  1,2,\ldots,n\right\}  $
satisfy $a_{j}\leq a_{\sigma\left(  i\right)  }$. In other words,
$a_{\sigma\left(  i\right)  }$ is an element $x$ of $\mathbb{Z}$ such that at
least $i$ elements $j\in\left\{  1,2,\ldots,n\right\}  $ satisfy $a_{j}\leq
x$. In other words,%
\[
a_{\sigma\left(  i\right)  }\in\left\{  x\in\mathbb{Z}\ \mid\ \text{at least
}i\text{ elements }j\in\left\{  1,2,\ldots,n\right\}  \text{ satisfy }%
a_{j}\leq x\right\}  =X.
\]

On the other hand, let $y$ be any element of $X$. We shall show that
$a_{\sigma\left(  i\right)  }\leq y$. Indeed, assume the contrary. Thus,
$a_{\sigma\left(  i\right)  }>y$. Hence,%
\begin{equation}
\left\vert \left\{  j\in\left\{  1,2,\ldots,n\right\}  \ \mid\ a_{j}\leq
y\right\}  \right\vert <i \label{pf.prop.sorting.bclaim.pf.5}%
\end{equation}
\footnote{\textit{Proof.} Let $p\in\left\{  j\in\left\{  1,2,\ldots,n\right\}
\ \mid\ a_{j}\leq y\right\}  $. Thus, $p$ is an element $j$ of $\left\{
1,2,\ldots,n\right\}  $ such that $a_{j}\leq y$. In other words, $p$ is an
element of $\left\{  1,2,\ldots,n\right\}  $ and satisfies $a_{p}\leq y$.
\par
The permutation $\sigma$ has an inverse $\sigma^{-1}$. Let $q=\sigma
^{-1}\left(  p\right)  $. Thus, $p=\sigma\left(  q\right)  $. Hence,
$a_{p}=a_{\sigma\left(  q\right)  }$, so that $a_{\sigma\left(  q\right)
}=a_{p}\leq y<a_{\sigma\left(  i\right)  }$ (since $a_{\sigma\left(  i\right)
}>y$).
\par
If we had $i\leq q$, then we would have $a_{\sigma\left(  i\right)  }\leq
a_{\sigma\left(  q\right)  }$ (by (\ref{pf.prop.sorting.bclaim.pf.1}), applied
to $u=i$ and $v=q$), which would contradict $a_{\sigma\left(  q\right)
}<a_{\sigma\left(  i\right)  }$. Hence, we cannot have $i\leq q$. Thus, we
must have $q<i$. Hence, $q\in\left\{  1,2,\ldots,i-1\right\}  $. Thus,
$p=\sigma\left(  \underbrace{q}_{\in\left\{  1,2,\ldots,i-1\right\}  }\right)
\in\sigma\left(  \left\{  1,2,\ldots,i-1\right\}  \right)  $.
\par
Let us now forget that we fixed $p$. We thus have shown that every
$p\in\left\{  j\in\left\{  1,2,\ldots,n\right\}  \ \mid\ a_{j}\leq y\right\}
$ satisfies $p\in\sigma\left(  \left\{  1,2,\ldots,i-1\right\}  \right)  $. In
other words,%
\[
\left\{  j\in\left\{  1,2,\ldots,n\right\}  \ \mid\ a_{j}\leq y\right\}
\subseteq\sigma\left(  \left\{  1,2,\ldots,i-1\right\}  \right)  .
\]
Hence,%
\begin{align*}
\left\vert \left\{  j\in\left\{  1,2,\ldots,n\right\}  \ \mid\ a_{j}\leq
y\right\}  \right\vert  &  \leq\left\vert \sigma\left(  \left\{
1,2,\ldots,i-1\right\}  \right)  \right\vert =\left\vert \left\{
1,2,\ldots,i-1\right\}  \right\vert \\
&  \ \ \ \ \ \ \ \ \ \ \left(  \text{since the map }\sigma\text{ is
injective}\right) \\
&  =i-1<i,
\end{align*}
qed.}.

But%
\[
y\in X=\left\{  x\in\mathbb{Z}\ \mid\ \text{at least }i\text{ elements }%
j\in\left\{  1,2,\ldots,n\right\}  \text{ satisfy }a_{j}\leq x\right\}  .
\]
In other words, $y$ is an element $x$ of $\mathbb{Z}$ such that at least $i$
elements $j\in\left\{  1,2,\ldots,n\right\}  $ satisfy $a_{j}\leq x$. In other
words, $y$ is an element of $\mathbb{Z}$, and at least $i$ elements
$j\in\left\{  1,2,\ldots,n\right\}  $ satisfy $a_{j}\leq y$. We have%
\[
\left\vert \left\{  j\in\left\{  1,2,\ldots,n\right\}  \ \mid\ a_{j}\leq
y\right\}  \right\vert \geq i
\]
(since at least $i$ elements $j\in\left\{  1,2,\ldots,n\right\}  $ satisfy
$a_{j}\leq y$). This contradicts (\ref{pf.prop.sorting.bclaim.pf.5}). This
contradiction proves that our assumption was wrong. Hence, $a_{\sigma\left(
i\right)  }\leq y$ is proven.

Now, let us forget that we fixed $y$. We thus have shown that every $y\in X$
satisfies $a_{\sigma\left(  i\right)  }\leq y$. Altogether, we thus have shown
the following two facts:

\begin{itemize}
\item We have $a_{\sigma\left(  i\right)  }\in X$.

\item Every $y\in X$ satisfies $a_{\sigma\left(  i\right)  }\leq y$.
\end{itemize}

Combining these two facts, we obtain $a_{\sigma\left(  i\right)  }=\min X$
(and, in particular, this shows that $\min X$ is well-defined). Since
\newline$X=\left\{  x\in\mathbb{Z}\ \mid\ \text{at least }i\text{ elements
}j\in\left\{  1,2,\ldots,n\right\}  \text{ satisfy }a_{j}\leq x\right\}  $,
this rewrites as follows:%
\[
a_{\sigma\left(  i\right)  }=\min\left\{  x\in\mathbb{Z}\ \mid\ \text{at least
}i\text{ elements }j\in\left\{  1,2,\ldots,n\right\}  \text{ satisfy }%
a_{j}\leq x\right\}  .
\]
Thus, (\ref{pf.prop.sorting.bclaim}) is proven.]

Now, the value $a_{\sigma\left(  i\right)  }$ depends only on $a_{1}%
,a_{2},\ldots,a_{n}$ and $i$ (but not on $\sigma$) (because
(\ref{pf.prop.sorting.bclaim}) gives a description of $a_{\sigma\left(
i\right)  }$ which involves $a_{1},a_{2},\ldots,a_{n}$ and $i$, but not
$\sigma$). Thus, Proposition \ref{prop.sorting} \textbf{(b)} is proven.

\textbf{(c)} The integers $a_{1},a_{2},\ldots,a_{n}$ are distinct. In other
words, if $u$ and $v$ are two distinct elements of $\left\{  1,2,\ldots
,n\right\}  $, then%
\begin{equation}
a_{u}\neq a_{v}. \label{pf.prop.sorting.c.distinct}%
\end{equation}

We have $\sigma\in S_{n}$. Thus, $\sigma$ is a permutation, thus a bijective
map, thus an injective map.

Now, we can see that:

\begin{itemize}
\item There is \textbf{at least one} permutation $\sigma\in S_{n}$ such that
$a_{\sigma\left(  1\right)  }<a_{\sigma\left(  2\right)  }<\cdots
<a_{\sigma\left(  n\right)  }$\ \ \ \ \footnote{\textit{Proof.} Proposition
\ref{prop.sorting} \textbf{(a)} shows that there exists a permutation
$\sigma\in S_{n}$ such that $a_{\sigma\left(  1\right)  }\leq a_{\sigma\left(
2\right)  }\leq\cdots\leq a_{\sigma\left(  n\right)  }$. Consider this
$\sigma$.
\par
Let $k\in\left\{  1,2,\ldots,n-1\right\}  $. Then, $a_{\sigma\left(  k\right)
}\leq a_{\sigma\left(  k+1\right)  }$ (since $a_{\sigma\left(  1\right)  }\leq
a_{\sigma\left(  2\right)  }\leq\cdots\leq a_{\sigma\left(  n\right)  }$). But
on the other hand, $k\neq k+1$. Hence, $\sigma\left(  k\right)  \neq
\sigma\left(  k+1\right)  $ (since the map $\sigma$ is injective). Hence,
$a_{\sigma\left(  k\right)  }\neq a_{\sigma\left(  k+1\right)  }$ (by
(\ref{pf.prop.sorting.c.distinct}), applied to $u=\sigma\left(  k\right)  $
and $v=\sigma\left(  k+1\right)  $). Combined with $a_{\sigma\left(  k\right)
}\leq a_{\sigma\left(  k+1\right)  }$, this yields $a_{\sigma\left(  k\right)
}<a_{\sigma\left(  k+1\right)  }$.
\par
Let us now forget that we fixed $k$. We thus have shown that $a_{\sigma\left(
k\right)  }<a_{\sigma\left(  k+1\right)  }$ for every $k\in\left\{
1,2,\ldots,n-1\right\}  $. In other words, $a_{\sigma\left(  1\right)
}<a_{\sigma\left(  2\right)  }<\cdots<a_{\sigma\left(  n\right)  }$.
\par
Let us now forget that we defined $\sigma$. We thus have constructed a
permutation $\sigma\in S_{n}$ such that $a_{\sigma\left(  1\right)
}<a_{\sigma\left(  2\right)  }<\cdots<a_{\sigma\left(  n\right)  }$. Hence,
there is \textbf{at least one} permutation $\sigma\in S_{n}$ such that
$a_{\sigma\left(  1\right)  }<a_{\sigma\left(  2\right)  }<\cdots
<a_{\sigma\left(  n\right)  }$. Qed.}.

\item There is \textbf{at most one} permutation $\sigma\in S_{n}$ such that
$a_{\sigma\left(  1\right)  }<a_{\sigma\left(  2\right)  }<\cdots
<a_{\sigma\left(  n\right)  }$\ \ \ \ \footnote{\textit{Proof.} Let
$\sigma_{1}$ and $\sigma_{2}$ be two permutations $\sigma\in S_{n}$ such that
$a_{\sigma\left(  1\right)  }<a_{\sigma\left(  2\right)  }<\cdots
<a_{\sigma\left(  n\right)  }$. We shall show that $\sigma_{1}=\sigma_{2}$.
\par
Fix $i\in\left\{  1,2,\ldots,n\right\}  $.
\par
We know that $\sigma_{1}$ is a permutation $\sigma\in S_{n}$ such that
$a_{\sigma\left(  1\right)  }<a_{\sigma\left(  2\right)  }<\cdots
<a_{\sigma\left(  n\right)  }$. In other words, $\sigma_{1}$ is a permutation
in $S_{n}$ such that $a_{\sigma_{1}\left(  1\right)  }<a_{\sigma_{1}\left(
2\right)  }<\cdots<a_{\sigma_{1}\left(  n\right)  }$. Thus,
(\ref{pf.prop.sorting.bclaim}) (applied to $\sigma=\sigma_{1}$) yields%
\begin{equation}
a_{\sigma_{1}\left(  i\right)  }=\min\left\{  x\in\mathbb{Z}\ \mid\ \text{at
least }i\text{ elements }j\in\left\{  1,2,\ldots,n\right\}  \text{ satisfy
}a_{j}\leq x\right\}  . \label{pf.prop.sorting.c.5}%
\end{equation}
The same argument (applied to $\sigma_{2}$ instead of $\sigma_{1}$) yields%
\[
a_{\sigma_{2}\left(  i\right)  }=\min\left\{  x\in\mathbb{Z}\ \mid\ \text{at
least }i\text{ elements }j\in\left\{  1,2,\ldots,n\right\}  \text{ satisfy
}a_{j}\leq x\right\}  .
\]
Comparing this with (\ref{pf.prop.sorting.c.5}), we obtain $a_{\sigma
_{1}\left(  i\right)  }=a_{\sigma_{2}\left(  i\right)  }$.
\par
Now, if we had $\sigma_{1}\left(  i\right)  \neq\sigma_{2}\left(  i\right)  $,
then we would have $a_{\sigma_{1}\left(  i\right)  }\neq a_{\sigma_{2}\left(
i\right)  }$ (by (\ref{pf.prop.sorting.c.distinct}), applied to $u=\sigma
_{1}\left(  i\right)  $ and $v=\sigma_{2}\left(  i\right)  $), which would
contradict $a_{\sigma_{1}\left(  i\right)  }=a_{\sigma_{2}\left(  i\right)  }%
$. Thus, we cannot have $\sigma_{1}\left(  i\right)  \neq\sigma_{2}\left(
i\right)  $. Hence, we must have $\sigma_{1}\left(  i\right)  =\sigma
_{2}\left(  i\right)  $.
\par
Let us now forget that we fixed $i$. We thus have shown that $\sigma
_{1}\left(  i\right)  =\sigma_{2}\left(  i\right)  $ for every $i\in\left\{
1,2,\ldots,n\right\}  $. In other words, $\sigma_{1}=\sigma_{2}$.
\par
Let us now forget that we fixed $\sigma_{1}$ and $\sigma_{2}$. We thus have
shown that $\sigma_{1}=\sigma_{2}$ whenever $\sigma_{1}$ and $\sigma_{2}$ are
two permutations $\sigma\in S_{n}$ such that $a_{\sigma\left(  1\right)
}<a_{\sigma\left(  2\right)  }<\cdots<a_{\sigma\left(  n\right)  }$. In other
words, there is \textbf{at most one} permutation $\sigma\in S_{n}$ such that
$a_{\sigma\left(  1\right)  }<a_{\sigma\left(  2\right)  }<\cdots
<a_{\sigma\left(  n\right)  }$. Qed.}.
\end{itemize}

Combining these two statements, we conclude that there is a \textbf{unique}
permutation $\sigma\in S_{n}$ such that $a_{\sigma\left(  1\right)
}<a_{\sigma\left(  2\right)  }<\cdots<a_{\sigma\left(  n\right)  }$. This
proves Proposition \ref{prop.sorting} \textbf{(c)}.
\end{proof}

\begin{proof}
[Proof of Lemma \ref{lem.cauchy-binet.EI}.]We have%
\[
\mathbf{E}=\left\{  \left(  k_{1},k_{2},\ldots,k_{n}\right)  \in\left[
m\right]  ^{n}\ \mid\ \text{the integers }k_{1},k_{2},\ldots,k_{n}\text{ are
distinct}\right\}
\]
and%
\[
\mathbf{I}=\left\{  \left(  k_{1},k_{2},\ldots,k_{n}\right)  \in\left[
m\right]  ^{n}\ \mid\ k_{1}<k_{2}<\cdots<k_{n}\right\}  .
\]

\begin{vershort}
Clearly, every element of $\mathbf{I}\times S_{n}$ can be written in the form
$\left(  \left(  g_{1},g_{2},\ldots,g_{n}\right)  ,\sigma\right)  $ for some
$\left(  g_{1},g_{2},\ldots,g_{n}\right)  \in\mathbf{I}$ and $\sigma\in S_{n}$.
\end{vershort}

\begin{verlong}
Every element of $\mathbf{I}\times S_{n}$ can be written in the form $\left(
\left(  g_{1},g_{2},\ldots,g_{n}\right)  ,\sigma\right)  $ for some $\left(
g_{1},g_{2},\ldots,g_{n}\right)  \in\mathbf{I}$ and $\sigma\in S_{n}%
$\ \ \ \ \footnote{\textit{Proof.} Let $\alpha\in\mathbf{I}\times S_{n}$.
Then, $\alpha$ can be written in the form $\alpha=\left(  \beta,\gamma\right)
$ for some $\beta\in\mathbf{I}$ and $\gamma\in S_{n}$ (since $\alpha
\in\mathbf{I}\times S_{n}$). Consider these $\beta$ and $\gamma$. We have
$\beta\in\mathbf{I}=\left\{  \left(  k_{1},k_{2},\ldots,k_{n}\right)
\in\left[  m\right]  ^{n}\ \mid\ k_{1}<k_{2}<\cdots<k_{n}\right\}
\subseteq\left[  m\right]  ^{n}$. Hence, $\beta$ can be written in the form
$\beta=\left(  k_{1},k_{2},\ldots,k_{n}\right)  $ for some $\left(
k_{1},k_{2},\ldots,k_{n}\right)  \in\left[  m\right]  ^{n}$. Consider this
$\left(  k_{1},k_{2},\ldots,k_{n}\right)  $. Clearly, $\left(  k_{1}%
,k_{2},\ldots,k_{n}\right)  =\beta\in\mathbf{I}$.
\par
But $\alpha=\left(  \underbrace{\beta}_{=\left(  k_{1},k_{2},\ldots
,k_{n}\right)  },\gamma\right)  =\left(  \left(  k_{1},k_{2},\ldots
,k_{n}\right)  ,\gamma\right)  $. Hence, $\alpha$ can be written in the form
$\left(  \left(  g_{1},g_{2},\ldots,g_{n}\right)  ,\sigma\right)  $ for some
$\left(  g_{1},g_{2},\ldots,g_{n}\right)  \in\mathbf{I}$ and $\sigma\in S_{n}$
(namely, for $\left(  g_{1},g_{2},\ldots,g_{n}\right)  =\left(  k_{1}%
,k_{2},\ldots,k_{n}\right)  $ and $\sigma=\gamma$).
\par
Now, let us forget that we fixed $\alpha$. We thus have proven that every
$\alpha\in\mathbf{I}\times S_{n}$ can be written in the form $\left(  \left(
g_{1},g_{2},\ldots,g_{n}\right)  ,\sigma\right)  $ for some $\left(
g_{1},g_{2},\ldots,g_{n}\right)  \in\mathbf{I}$ and $\sigma\in S_{n}$. In
other words, every element of $\mathbf{I}\times S_{n}$ can be written in the
form $\left(  \left(  g_{1},g_{2},\ldots,g_{n}\right)  ,\sigma\right)  $ for
some $\left(  g_{1},g_{2},\ldots,g_{n}\right)  \in\mathbf{I}$ and $\sigma\in
S_{n}$. Qed.}.
\end{verlong}

For every $\left(  \left(  g_{1},g_{2},\ldots,g_{n}\right)  ,\sigma\right)
\in\mathbf{I}\times S_{n}$, we have $\left(  g_{\sigma\left(  1\right)
},g_{\sigma\left(  2\right)  },\ldots,g_{\sigma\left(  n\right)  }\right)
\in\mathbf{E}$\ \ \ \ \footnote{\textit{Proof.} Let $\left(  \left(
g_{1},g_{2},\ldots,g_{n}\right)  ,\sigma\right)  \in\mathbf{I}\times S_{n}$.
Thus, $\left(  g_{1},g_{2},\ldots,g_{n}\right)  \in\mathbf{I}$ and $\sigma\in
S_{n}$.
\par
We have $\sigma\in S_{n}$. Thus, $\sigma$ is a permutation, and thus a
bijective map, hence an injective map.
\par
We have
\[
\left(  g_{1},g_{2},\ldots,g_{n}\right)  \in\mathbf{I}=\left\{  \left(
k_{1},k_{2},\ldots,k_{n}\right)  \in\left[  m\right]  ^{n}\ \mid\ k_{1}%
<k_{2}<\cdots<k_{n}\right\}  .
\]
In other words, $\left(  g_{1},g_{2},\ldots,g_{n}\right)  $ is an element
$\left(  k_{1},k_{2},\ldots,k_{n}\right)  \in\left[  m\right]  ^{n}$
satisfying $k_{1}<k_{2}<\cdots<k_{n}$. In other words, $\left(  g_{1}%
,g_{2},\ldots,g_{n}\right)  $ is an element of $\left[  m\right]  ^{n}$ and
satisfies $g_{1}<g_{2}<\cdots<g_{n}$. Hence, the integers $g_{1},g_{2}%
,\ldots,g_{n}$ are distinct (since $g_{1}<g_{2}<\cdots<g_{n}$). In other
words, any two distinct elements $u$ and $v$ of $\left[  n\right]  $ satisfy%
\begin{equation}
g_{u}\neq g_{v}. \label{pf.lem.cauchy-binet.EI.1.pf.1}%
\end{equation}
\par
Now, let $u$ and $v$ be two distinct elements of $\left[  n\right]  $. Thus,
$u\neq v$ (since $u$ and $v$ are distinct), so that $\sigma\left(  u\right)
\neq\sigma\left(  v\right)  $ (since $\sigma$ is injective). Hence,
$g_{\sigma\left(  u\right)  }\neq g_{\sigma\left(  v\right)  }$ (by
(\ref{pf.lem.cauchy-binet.EI.1.pf.1}), applied to $\sigma\left(  u\right)  $
and $\sigma\left(  v\right)  $ instead of $u$ and $v$).
\par
Let us now forget that we fixed $u$ and $v$. We thus have shown that any two
distinct elements $u$ and $v$ of $\left[  n\right]  $ satisfy $g_{\sigma
\left(  u\right)  }\neq g_{\sigma\left(  v\right)  }$. In other words, the
integers $g_{\sigma\left(  1\right)  },g_{\sigma\left(  2\right)  }%
,\ldots,g_{\sigma\left(  n\right)  }$ are distinct. Hence, $\left(
g_{\sigma\left(  1\right)  },g_{\sigma\left(  2\right)  },\ldots
,g_{\sigma\left(  n\right)  }\right)  $ is an element $\left(  k_{1}%
,k_{2},\ldots,k_{n}\right)  \in\left[  m\right]  ^{n}$ such that the integers
$k_{1},k_{2},\ldots,k_{n}$ are distinct. In other words,%
\[
\left(  g_{\sigma\left(  1\right)  },g_{\sigma\left(  2\right)  }%
,\ldots,g_{\sigma\left(  n\right)  }\right)  \in\left\{  \left(  k_{1}%
,k_{2},\ldots,k_{n}\right)  \in\left[  m\right]  ^{n}\ \mid\ \text{the
integers }k_{1},k_{2},\ldots,k_{n}\text{ are distinct}\right\}  =\mathbf{E},
\]
qed.}. Hence, we can define a map $\Phi:\mathbf{I}\times S_{n}\rightarrow
\mathbf{E}$ by setting%
\begin{equation}
\left(
\begin{array}
[c]{c}%
\Phi\left(  \left(  g_{1},g_{2},\ldots,g_{n}\right)  ,\sigma\right)  =\left(
g_{\sigma\left(  1\right)  },g_{\sigma\left(  2\right)  },\ldots
,g_{\sigma\left(  n\right)  }\right) \\
\ \ \ \ \ \ \ \ \ \ \text{for every }\left(  \left(  g_{1},g_{2},\ldots
,g_{n}\right)  ,\sigma\right)  \in\mathbf{I}\times S_{n}%
\end{array}
\right)  \label{pf.lem.cauchy-binet.EI.defPhi}%
\end{equation}
(since every element of $\mathbf{I}\times S_{n}$ can be written in the form
$\left(  \left(  g_{1},g_{2},\ldots,g_{n}\right)  ,\sigma\right)  $ for some
$\left(  g_{1},g_{2},\ldots,g_{n}\right)  \in\mathbf{I}$ and $\sigma\in S_{n}%
$). Consider this map $\Phi$.

The map $\Phi$ is the map
\begin{align*}
\mathbf{I}\times S_{n}  &  \rightarrow\mathbf{E},\\
\left(  \left(  g_{1},g_{2},\ldots,g_{n}\right)  ,\sigma\right)   &
\mapsto\left(  g_{\sigma\left(  1\right)  },g_{\sigma\left(  2\right)
},\ldots,g_{\sigma\left(  n\right)  }\right)
\end{align*}
(since $\Phi$ is defined by (\ref{pf.lem.cauchy-binet.EI.defPhi})). In
particular, the latter map is well-defined.

The map $\Phi$ is injective\footnote{\textit{Proof.} Let $\alpha$ and $\beta$
be two elements of $\mathbf{I}\times S_{n}$ such that $\Phi\left(
\alpha\right)  =\Phi\left(  \beta\right)  $. We shall show that $\alpha=\beta
$.
\par
Let $\gamma$ be the element $\Phi\left(  \alpha\right)  =\Phi\left(
\beta\right)  $ of $\mathbf{E}$. Thus, $\gamma=\Phi\left(  \alpha\right)
=\Phi\left(  \beta\right)  $.
\par
Write $\alpha\in\mathbf{I}\times S_{n}$ in the form $\alpha=\left(  \left(
g_{1},g_{2},\ldots,g_{n}\right)  ,\pi\right)  $ for some $\left(  g_{1}%
,g_{2},\ldots,g_{n}\right)  \in\mathbf{I}$ and $\pi\in S_{n}$.
\par
Write $\beta\in\mathbf{I}\times S_{n}$ in the form $\beta=\left(  \left(
h_{1},h_{2},\ldots,h_{n}\right)  ,\tau\right)  $ for some $\left(  h_{1}%
,h_{2},\ldots,h_{n}\right)  \in\mathbf{I}$ and $\tau\in S_{n}$.
\par
We have%
\[
\left(  g_{1},g_{2},\ldots,g_{n}\right)  \in\mathbf{I}=\left\{  \left(
k_{1},k_{2},\ldots,k_{n}\right)  \in\left[  m\right]  ^{n}\ \mid\ k_{1}%
<k_{2}<\cdots<k_{n}\right\}  .
\]
In other words, $\left(  g_{1},g_{2},\ldots,g_{n}\right)  $ is an element
$\left(  k_{1},k_{2},\ldots,k_{n}\right)  $ of $\left[  m\right]  ^{n}$
satisfying $k_{1}<k_{2}<\cdots<k_{n}$. In other words, $\left(  g_{1}%
,g_{2},\ldots,g_{n}\right)  $ is an element of $\left[  m\right]  ^{n}$ and
satisfies $g_{1}<g_{2}<\cdots<g_{n}$.
\par
We have%
\[
\left(  h_{1},h_{2},\ldots,h_{n}\right)  \in\mathbf{I}=\left\{  \left(
k_{1},k_{2},\ldots,k_{n}\right)  \in\left[  m\right]  ^{n}\ \mid\ k_{1}%
<k_{2}<\cdots<k_{n}\right\}  .
\]
In other words, $\left(  h_{1},h_{2},\ldots,h_{n}\right)  $ is an element
$\left(  k_{1},k_{2},\ldots,k_{n}\right)  $ of $\left[  m\right]  ^{n}$
satisfying $k_{1}<k_{2}<\cdots<k_{n}$. In other words, $\left(  h_{1}%
,h_{2},\ldots,h_{n}\right)  $ is an element of $\left[  m\right]  ^{n}$ and
satisfies $h_{1}<h_{2}<\cdots<h_{n}$.
\par
We have $\gamma\in\mathbf{E}=\left\{  \left(  k_{1},k_{2},\ldots,k_{n}\right)
\in\left[  m\right]  ^{n}\ \mid\ \text{the integers }k_{1},k_{2},\ldots
,k_{n}\text{ are distinct}\right\}  $. In other words, we can write $\gamma$
in the form $\gamma=\left(  k_{1},k_{2},\ldots,k_{n}\right)  $ for some
$\left(  k_{1},k_{2},\ldots,k_{n}\right)  \in\left[  m\right]  ^{n}$ such that
the integers $k_{1},k_{2},\ldots,k_{n}$ are distinct. Consider this $\left(
k_{1},k_{2},\ldots,k_{n}\right)  $.
\par
Applying the map $\Phi$ to both sides of the equality $\alpha=\left(  \left(
g_{1},g_{2},\ldots,g_{n}\right)  ,\pi\right)  $, we obtain%
\[
\Phi\left(  \alpha\right)  =\Phi\left(  \left(  g_{1},g_{2},\ldots
,g_{n}\right)  ,\pi\right)  =\left(  g_{\pi\left(  1\right)  },g_{\pi\left(
2\right)  },\ldots,g_{\pi\left(  n\right)  }\right)
\ \ \ \ \ \ \ \ \ \ \left(  \text{by the definition of }\Phi\right)  .
\]
Thus, $\left(  g_{\pi\left(  1\right)  },g_{\pi\left(  2\right)  }%
,\ldots,g_{\pi\left(  n\right)  }\right)  =\Phi\left(  \alpha\right)
=\gamma=\left(  k_{1},k_{2},\ldots,k_{n}\right)  $. In other words, every
$i\in\left\{  1,2,\ldots,n\right\}  $ satisfies
\begin{equation}
g_{\pi\left(  i\right)  }=k_{i}. \label{pf.lem.cauchy-binet.EI.2.pf.1}%
\end{equation}
Thus, every $j\in\left\{  1,2,\ldots,n\right\}  $ satisfies
\begin{align}
g_{j}  &  =g_{\pi\left(  \pi^{-1}\left(  j\right)  \right)  }%
\ \ \ \ \ \ \ \ \ \ \left(  \text{since }j=\pi\left(  \pi^{-1}\left(
j\right)  \right)  \right) \nonumber\\
&  =k_{\pi^{-1}\left(  j\right)  }\ \ \ \ \ \ \ \ \ \ \left(  \text{by
(\ref{pf.lem.cauchy-binet.EI.2.pf.1}), applied to }i=\pi^{-1}\left(  j\right)
\right)  . \label{pf.lem.cauchy-binet.EI.2.pf.1a}%
\end{align}
\par
We have $g_{1}<g_{2}<\cdots<g_{n}$. This rewrites as
\[
k_{\pi^{-1}\left(  1\right)  }<k_{\pi^{-1}\left(  2\right)  }<\cdots
<k_{\pi^{-1}\left(  n\right)  }%
\]
(because every $j\in\left\{  1,2,\ldots,n\right\}  $ satisfies $g_{j}%
=k_{\pi^{-1}\left(  j\right)  }$).
\par
Applying the map $\Phi$ to both sides of the equality $\beta=\left(  \left(
h_{1},h_{2},\ldots,h_{n}\right)  ,\tau\right)  $, we obtain%
\[
\Phi\left(  \beta\right)  =\Phi\left(  \left(  h_{1},h_{2},\ldots
,h_{n}\right)  ,\tau\right)  =\left(  h_{\tau\left(  1\right)  }%
,h_{\tau\left(  2\right)  },\ldots,h_{\tau\left(  n\right)  }\right)
\ \ \ \ \ \ \ \ \ \ \left(  \text{by the definition of }\Phi\right)  .
\]
Thus, $\left(  h_{\tau\left(  1\right)  },h_{\tau\left(  2\right)  }%
,\ldots,h_{\tau\left(  n\right)  }\right)  =\Phi\left(  \beta\right)
=\gamma=\left(  k_{1},k_{2},\ldots,k_{n}\right)  $. In other words, every
$i\in\left\{  1,2,\ldots,n\right\}  $ satisfies
\begin{equation}
h_{\tau\left(  i\right)  }=k_{i}. \label{pf.lem.cauchy-binet.EI.2.pf.2}%
\end{equation}
Thus, every $j\in\left\{  1,2,\ldots,n\right\}  $ satisfies
\begin{align}
h_{j}  &  =h_{\tau\left(  \tau^{-1}\left(  j\right)  \right)  }%
\ \ \ \ \ \ \ \ \ \ \left(  \text{since }j=\tau\left(  \tau^{-1}\left(
j\right)  \right)  \right) \nonumber\\
&  =k_{\tau^{-1}\left(  j\right)  }\ \ \ \ \ \ \ \ \ \ \left(  \text{by
(\ref{pf.lem.cauchy-binet.EI.2.pf.2}), applied to }i=\tau^{-1}\left(
j\right)  \right)  . \label{pf.lem.cauchy-binet.EI.2.pf.2a}%
\end{align}
\par
We have $h_{1}<h_{2}<\cdots<h_{n}$. This rewrites as%
\[
k_{\tau^{-1}\left(  1\right)  }<k_{\tau^{-1}\left(  2\right)  }<\cdots
<k_{\tau^{-1}\left(  n\right)  }%
\]
(because every $j\in\left\{  1,2,\ldots,n\right\}  $ satisfies $h_{j}%
=k_{\tau^{-1}\left(  j\right)  }$).
\par
Proposition \ref{prop.sorting} \textbf{(c)} (applied to $a_{i}=k_{i}$) yields
that there is a \textbf{unique} permutation $\sigma\in S_{n}$ such that
$k_{\sigma\left(  1\right)  }<k_{\sigma\left(  2\right)  }<\cdots
<k_{\sigma\left(  n\right)  }$ (since the integers $k_{1},k_{2},\ldots,k_{n}$
are distinct). In particular, there exists \textbf{at most one} such
permutation. In other words, if $\sigma_{1}$ and $\sigma_{2}$ are two
permutations $\sigma\in S_{n}$ satisfying $k_{\sigma\left(  1\right)
}<k_{\sigma\left(  2\right)  }<\cdots<k_{\sigma\left(  n\right)  }$, then%
\begin{equation}
\sigma_{1}=\sigma_{2}. \label{pf.lem.cauchy-binet.EI.2.pf.4}%
\end{equation}
\par
Now, $\pi^{-1}$ is a permutation $\sigma\in S_{n}$ satisfying $k_{\sigma
\left(  1\right)  }<k_{\sigma\left(  2\right)  }<\cdots<k_{\sigma\left(
n\right)  }$ (since $k_{\pi^{-1}\left(  1\right)  }<k_{\pi^{-1}\left(
2\right)  }<\cdots<k_{\pi^{-1}\left(  n\right)  }$). Also, $\tau^{-1}$ is a
permutation $\sigma\in S_{n}$ satisfying $k_{\sigma\left(  1\right)
}<k_{\sigma\left(  2\right)  }<\cdots<k_{\sigma\left(  n\right)  }$ (since
$k_{\tau^{-1}\left(  1\right)  }<k_{\tau^{-1}\left(  2\right)  }%
<\cdots<k_{\tau^{-1}\left(  n\right)  }$). Hence, we can apply
(\ref{pf.lem.cauchy-binet.EI.2.pf.4}) to $\sigma_{1}=\pi^{-1}$ and $\sigma
_{2}=\tau^{-1}$. As a result, we obtain $\pi^{-1}=\tau^{-1}$. Hence, $\pi
=\tau$. Now, every $j\in\left\{  1,2,\ldots,n\right\}  $ satisfies%
\begin{align*}
g_{j}  &  =k_{\pi^{-1}\left(  j\right)  }\ \ \ \ \ \ \ \ \ \ \left(  \text{by
(\ref{pf.lem.cauchy-binet.EI.2.pf.1a})}\right) \\
&  =k_{\tau^{-1}\left(  j\right)  }\ \ \ \ \ \ \ \ \ \ \left(  \text{since
}\pi=\tau\right) \\
&  =h_{j}\ \ \ \ \ \ \ \ \ \ \left(  \text{by
(\ref{pf.lem.cauchy-binet.EI.2.pf.2a})}\right)  .
\end{align*}
In other words, $\left(  g_{1},g_{2},\ldots,g_{n}\right)  =\left(  h_{1}%
,h_{2},\ldots,h_{n}\right)  $. Thus,
\[
\alpha=\left(  \underbrace{\left(  g_{1},g_{2},\ldots,g_{n}\right)
}_{=\left(  h_{1},h_{2},\ldots,h_{n}\right)  },\underbrace{\pi}_{=\tau
}\right)  =\left(  \left(  h_{1},h_{2},\ldots,h_{n}\right)  ,\tau\right)
=\beta.
\]
\par
Now, let us forget that we fixed $\alpha$ and $\beta$. We thus have shown that
if $\alpha$ and $\beta$ are two elements of $\mathbf{I}\times S_{n}$ such that
$\Phi\left(  \alpha\right)  =\Phi\left(  \beta\right)  $, then $\alpha=\beta$.
In other words, the map $\Phi$ is injective, qed.} and
surjective\footnote{\textit{Proof.} Let $\gamma\in\mathbf{E}$. We shall prove
that $\gamma\in\Phi\left(  \mathbf{I}\times S_{n}\right)  $.
\par
We have $\gamma\in\mathbf{E}=\left\{  \left(  k_{1},k_{2},\ldots,k_{n}\right)
\in\left[  m\right]  ^{n}\ \mid\ \text{the integers }k_{1},k_{2},\ldots
,k_{n}\text{ are distinct}\right\}  $. In other words, we can write $\gamma$
in the form $\gamma=\left(  k_{1},k_{2},\ldots,k_{n}\right)  $ for some
$\left(  k_{1},k_{2},\ldots,k_{n}\right)  \in\left[  m\right]  ^{n}$ such that
the integers $k_{1},k_{2},\ldots,k_{n}$ are distinct. Let us denote this
$\left(  k_{1},k_{2},\ldots,k_{n}\right)  $ by $\left(  a_{1},a_{2}%
,\ldots,a_{n}\right)  $. Thus, $\left(  a_{1},a_{2},\ldots,a_{n}\right)  $ is
an element of $\left[  m\right]  ^{n}$ such that the integers $a_{1}%
,a_{2},\ldots,a_{n}$ are distinct, and we have $\gamma=\left(  a_{1}%
,a_{2},\ldots,a_{n}\right)  $.
\par
Proposition \ref{prop.sorting} \textbf{(c)} yields that there is a
\textbf{unique} permutation $\sigma\in S_{n}$ such that $a_{\sigma\left(
1\right)  }<a_{\sigma\left(  2\right)  }<\cdots<a_{\sigma\left(  n\right)  }$
(since the integers $a_{1},a_{2},\ldots,a_{n}$ are distinct). In particular,
there exists \textbf{at least one} such permutation. Consider such a
permutation $\sigma$. Thus, $\sigma\in S_{n}$ and $a_{\sigma\left(  1\right)
}<a_{\sigma\left(  2\right)  }<\cdots<a_{\sigma\left(  n\right)  }$.
\par
Now, $\left(  a_{\sigma\left(  1\right)  },a_{\sigma\left(  2\right)  }%
,\ldots,a_{\sigma\left(  n\right)  }\right)  $ is an element of $\left[
m\right]  ^{n}$ satisfying $a_{\sigma\left(  1\right)  }<a_{\sigma\left(
2\right)  }<\cdots<a_{\sigma\left(  n\right)  }$. In other words, $\left(
a_{\sigma\left(  1\right)  },a_{\sigma\left(  2\right)  },\ldots
,a_{\sigma\left(  n\right)  }\right)  $ is an element $\left(  k_{1}%
,k_{2},\ldots,k_{n}\right)  $ of $\left[  m\right]  ^{n}$ satisfying
$k_{1}<k_{2}<\cdots<k_{n}$. In other words,
\[
\left(  a_{\sigma\left(  1\right)  },a_{\sigma\left(  2\right)  }%
,\ldots,a_{\sigma\left(  n\right)  }\right)  \in\left\{  \left(  k_{1}%
,k_{2},\ldots,k_{n}\right)  \in\left[  m\right]  ^{n}\ \mid\ k_{1}%
<k_{2}<\cdots<k_{n}\right\}  =\mathbf{I}.
\]
Hence, $\Phi\left(  \left(  a_{\sigma\left(  1\right)  },a_{\sigma\left(
2\right)  },\ldots,a_{\sigma\left(  n\right)  }\right)  ,\sigma^{-1}\right)  $
is well-defined. The definition of $\Phi\left(  \left(  a_{\sigma\left(
1\right)  },a_{\sigma\left(  2\right)  },\ldots,a_{\sigma\left(  n\right)
}\right)  ,\sigma^{-1}\right)  $ yields%
\begin{align*}
&  \Phi\left(  \left(  a_{\sigma\left(  1\right)  },a_{\sigma\left(  2\right)
},\ldots,a_{\sigma\left(  n\right)  }\right)  ,\sigma^{-1}\right) \\
&  =\left(  a_{\sigma\left(  \sigma^{-1}\left(  1\right)  \right)  }%
,a_{\sigma\left(  \sigma^{-1}\left(  2\right)  \right)  },\ldots
,a_{\sigma\left(  \sigma^{-1}\left(  n\right)  \right)  }\right)  =\left(
a_{1},a_{2},\ldots,a_{n}\right) \\
&  \ \ \ \ \ \ \ \ \ \ \left(  \text{since }a_{\sigma\left(  \sigma
^{-1}\left(  i\right)  \right)  }=a_{i}\text{ for every }i\in\left\{
1,2,\ldots,n\right\}  \right) \\
&  =\gamma.
\end{align*}
Thus, $\gamma=\Phi\left(  \underbrace{\left(  a_{\sigma\left(  1\right)
},a_{\sigma\left(  2\right)  },\ldots,a_{\sigma\left(  n\right)  }\right)
}_{\in\mathbf{I}},\underbrace{\sigma^{-1}}_{\in S_{n}}\right)  \in\Phi\left(
\mathbf{I}\times S_{n}\right)  $.
\par
Now, let us forget that we fixed $\gamma$. We thus have shown that $\gamma
\in\Phi\left(  \mathbf{I}\times S_{n}\right)  $ for every $\gamma\in
\mathbf{E}$. In other words, $\mathbf{E}\subseteq\Phi\left(  \mathbf{I}\times
S_{n}\right)  $. In other words, the map $\Phi$ is surjective. Qed.}. Hence,
the map $\Phi$ is bijective. In other words, the map $\Phi$ is a bijection. In
other words, the map
\begin{align*}
\mathbf{I}\times S_{n}  &  \rightarrow\mathbf{E},\\
\left(  \left(  g_{1},g_{2},\ldots,g_{n}\right)  ,\sigma\right)   &
\mapsto\left(  g_{\sigma\left(  1\right)  },g_{\sigma\left(  2\right)
},\ldots,g_{\sigma\left(  n\right)  }\right)
\end{align*}
is a bijection (because the map $\Phi$ is the map
\begin{align*}
\mathbf{I}\times S_{n}  &  \rightarrow\mathbf{E},\\
\left(  \left(  g_{1},g_{2},\ldots,g_{n}\right)  ,\sigma\right)   &
\mapsto\left(  g_{\sigma\left(  1\right)  },g_{\sigma\left(  2\right)
},\ldots,g_{\sigma\left(  n\right)  }\right)
\end{align*}
). This concludes the proof of Lemma \ref{lem.cauchy-binet.EI}.
\end{proof}

\subsection{Solution to Exercise \ref{exe.sorting.nmu}}

Before we start solving Exercise \ref{exe.sorting.nmu}, let us isolate a
useful fact that was proven in our above proof of Proposition
\ref{prop.sorting} \textbf{(b)}:

\begin{lemma}
\label{lem.sorting.nmu.uniprop}Let $n\in\mathbb{N}$. Let $a_{1},a_{2}%
,\ldots,a_{n}$ be $n$ integers. Let $\sigma\in S_{n}$ be such that
$a_{\sigma\left(  1\right)  }\leq a_{\sigma\left(  2\right)  }\leq\cdots\leq
a_{\sigma\left(  n\right)  }$. Let $i\in\left\{  1,2,\ldots,n\right\}  $.
Then,%
\begin{align}
&  a_{\sigma\left(  i\right)  }\nonumber\\
&  =\min\left\{  x\in\mathbb{Z}\ \mid\ \text{at least }i\text{ elements }%
j\in\left\{  1,2,\ldots,n\right\}  \text{ satisfy }a_{j}\leq x\right\}
\ \ \ \ \ \ \label{eq.lem.sorting.nmu.uniprop.claim}%
\end{align}
(and, in particular, the right hand side of
(\ref{eq.lem.sorting.nmu.uniprop.claim}) is well-defined).
\end{lemma}

\begin{proof}
[Proof of Lemma \ref{lem.sorting.nmu.uniprop}.]The equality
(\ref{eq.lem.sorting.nmu.uniprop.claim}) is precisely the equality
(\ref{pf.prop.sorting.bclaim}) that was proven during our above proof of
Proposition \ref{prop.sorting} \textbf{(b)}. Hence, we do not need to prove it
again. Thus, Lemma \ref{lem.sorting.nmu.uniprop} is proven.
\end{proof}

\begin{proof}
[First solution to Exercise \ref{exe.sorting.nmu}.]Let $i\in\left\{
1,2,\ldots,m\right\}  $. We must prove that $a_{\sigma\left(  i\right)  }\leq
b_{\tau\left(  i\right)  }$.

We have $n\geq m$, so that $m\leq n$. Now, $i\in\left\{  1,2,\ldots,m\right\}
\subseteq\left\{  1,2,\ldots,n\right\}  $ (since $m\leq n$). Hence,
$\sigma\left(  i\right)  $ is well-defined (since $\sigma$ is a map $\left\{
1,2,\ldots,n\right\}  \rightarrow\left\{  1,2,\ldots,n\right\}  $ (since
$\sigma\in S_{n}$)).

Define a subset $X$ of $\mathbb{Z}$ by%
\begin{equation}
X=\left\{  x\in\mathbb{Z}\ \mid\ \text{at least }i\text{ elements }%
j\in\left\{  1,2,\ldots,n\right\}  \text{ satisfy }a_{j}\leq x\right\}  .
\label{sol.sorting.nmu.X=}%
\end{equation}
Define a subset $Y$ of $\mathbb{Z}$ by%
\begin{equation}
Y=\left\{  x\in\mathbb{Z}\ \mid\ \text{at least }i\text{ elements }%
j\in\left\{  1,2,\ldots,m\right\}  \text{ satisfy }b_{j}\leq x\right\}  .
\label{sol.sorting.nmu.Y=}%
\end{equation}

Lemma \ref{lem.sorting.nmu.uniprop} yields that
\begin{align}
&  a_{\sigma\left(  i\right)  }\nonumber\\
&  =\min\left\{  x\in\mathbb{Z}\ \mid\ \text{at least }i\text{ elements }%
j\in\left\{  1,2,\ldots,n\right\}  \text{ satisfy }a_{j}\leq x\right\}
\ \ \ \ \ \ \label{sol.sorting.nmu.uniprop1}%
\end{align}
(and, in particular, the right hand side of (\ref{sol.sorting.nmu.uniprop1})
is well-defined). Thus,%
\begin{align}
a_{\sigma\left(  i\right)  }  &  =\min\underbrace{\left\{  x\in\mathbb{Z}%
\ \mid\ \text{at least }i\text{ elements }j\in\left\{  1,2,\ldots,n\right\}
\text{ satisfy }a_{j}\leq x\right\}  }_{\substack{=X\\\text{(by
(\ref{sol.sorting.nmu.X=}))}}}\nonumber\\
&  =\min X. \label{sol.sorting.nmu.uniprop1'}%
\end{align}

But $\tau\in S_{m}$ satisfies $b_{\tau\left(  1\right)  }\leq b_{\tau\left(
2\right)  }\leq\cdots\leq b_{\tau\left(  m\right)  }$, and we have
$i\in\left\{  1,2,\ldots,m\right\}  $. Thus, Lemma
\ref{lem.sorting.nmu.uniprop} (applied to $m$, $b_{k}$ and $\tau$ instead of
$n$, $a_{k}$ and $\sigma$) yields that%
\begin{equation}
b_{\tau\left(  i\right)  }=\min\left\{  x\in\mathbb{Z}\ \mid\ \text{at least
}i\text{ elements }j\in\left\{  1,2,\ldots,m\right\}  \text{ satisfy }%
b_{j}\leq x\right\}  \label{sol.sorting.nmu.uniprop2}%
\end{equation}
(and, in particular, the right hand side of (\ref{sol.sorting.nmu.uniprop2})
is well-defined). Thus,%
\begin{align}
b_{\tau\left(  i\right)  }  &  =\min\underbrace{\left\{  x\in\mathbb{Z}%
\ \mid\ \text{at least }i\text{ elements }j\in\left\{  1,2,\ldots,m\right\}
\text{ satisfy }b_{j}\leq x\right\}  }_{\substack{=Y\\\text{(by
(\ref{sol.sorting.nmu.Y=}))}}}\nonumber\\
&  =\min Y. \label{sol.sorting.nmu.uniprop2'}%
\end{align}

\begin{vershort}
But clearly, $\min Y\in Y$. Hence,
\[
b_{\tau\left(  i\right)  }=\min Y\in Y=\left\{  x\in\mathbb{Z}\ \mid\ \text{at
least }i\text{ elements }j\in\left\{  1,2,\ldots,m\right\}  \text{ satisfy
}b_{j}\leq x\right\}  .
\]
In other words, $b_{\tau\left(  i\right)  }$ is an element of $\mathbb{Z}$
such that at least $i$ elements $j\in\left\{  1,2,\ldots,m\right\}  $ satisfy
$b_{j}\leq b_{\tau\left(  i\right)  }$.
\end{vershort}

\begin{verlong}
But if $S$ is a subset of $\mathbb{Z}$ for which $\min S$ is well-defined,
then $\min S\in S$ (since the minimum of a set must lie inside the set).
Applying this to $S=Y$, we obtain $\min Y\in Y$. Now,
(\ref{sol.sorting.nmu.uniprop2'}) becomes%
\[
b_{\tau\left(  i\right)  }=\min Y\in Y=\left\{  x\in\mathbb{Z}\ \mid\ \text{at
least }i\text{ elements }j\in\left\{  1,2,\ldots,m\right\}  \text{ satisfy
}b_{j}\leq x\right\}  .
\]
In other words, $b_{\tau\left(  i\right)  }$ is an element $x$ of $\mathbb{Z}$
such that at least $i$ elements $j\in\left\{  1,2,\ldots,m\right\}  $ satisfy
$b_{j}\leq x$. In other words, $b_{\tau\left(  i\right)  }$ is an element of
$\mathbb{Z}$ such that at least $i$ elements $j\in\left\{  1,2,\ldots
,m\right\}  $ satisfy $b_{j}\leq b_{\tau\left(  i\right)  }$.
\end{verlong}

\begin{vershort}
We now know that at least $i$ elements $j\in\left\{  1,2,\ldots,m\right\}  $
satisfy $b_{j}\leq b_{\tau\left(  i\right)  }$. Each of these $i$ elements $j$
must also be an element of $\left\{  1,2,\ldots,n\right\}  $ (since
$j\in\left\{  1,2,\ldots,m\right\}  \subseteq\left\{  1,2,\ldots,n\right\}  $)
which satisfies $a_{j}\leq b_{\tau\left(  i\right)  }$ (because
(\ref{eq.exe.sorting.nmu.ass}) (applied to $j$ instead of $i$) shows that
$a_{j}\leq b_{j}\leq b_{\tau\left(  i\right)  }$). Thus, at least $i$ elements
$j\in\left\{  1,2,\ldots,n\right\}  $ satisfy $a_{j}\leq b_{\tau\left(
i\right)  }$.
\end{vershort}

\begin{verlong}
We now know that at least $i$ elements $j\in\left\{  1,2,\ldots,m\right\}  $
satisfy $b_{j}\leq b_{\tau\left(  i\right)  }$. In other words,%
\[
\left\vert \left\{  j\in\left\{  1,2,\ldots,m\right\}  \ \mid\ b_{j}\leq
b_{\tau\left(  i\right)  }\right\}  \right\vert \geq i.
\]

But $\left\{  j\in\left\{  1,2,\ldots,m\right\}  \ \mid\ b_{j}\leq
b_{\tau\left(  i\right)  }\right\}  \subseteq\left\{  j\in\left\{
1,2,\ldots,n\right\}  \ \mid\ a_{j}\leq b_{\tau\left(  i\right)  }\right\}
$\ \ \ \ \footnote{\textit{Proof.} Let $g\in\left\{  j\in\left\{
1,2,\ldots,m\right\}  \ \mid\ b_{j}\leq b_{\tau\left(  i\right)  }\right\}  $.
We shall prove that $g\in\left\{  j\in\left\{  1,2,\ldots,m\right\}
\ \mid\ a_{j}\leq b_{\tau\left(  i\right)  }\right\}  $.
\par
We have $g\in\left\{  j\in\left\{  1,2,\ldots,m\right\}  \ \mid\ b_{j}\leq
b_{\tau\left(  i\right)  }\right\}  $. In other words, $g$ is an element $j$
of $\left\{  1,2,\ldots,m\right\}  $ satisfying $b_{j}\leq b_{\tau\left(
i\right)  }$. In other words, $g$ is an element of $\left\{  1,2,\ldots
,m\right\}  $ satisfying $b_{g}\leq b_{\tau\left(  i\right)  }$. Now,
$g\in\left\{  1,2,\ldots,m\right\}  \subseteq\left\{  1,2,\ldots,n\right\}  $.
Also, (\ref{eq.exe.sorting.nmu.ass}) (applied to $g$ instead of $i$) shows
that $a_{g}\leq b_{g}\leq b_{\tau\left(  i\right)  }$. Hence, $g$ is an
element of $\left\{  1,2,\ldots,n\right\}  $ satisfying $a_{g}\leq
b_{\tau\left(  i\right)  }$. In other words, $g$ is an element $j$ of
$\left\{  1,2,\ldots,n\right\}  $ satisfying $a_{j}\leq b_{\tau\left(
i\right)  }$. In other words, $g\in\left\{  j\in\left\{  1,2,\ldots,n\right\}
\ \mid\ a_{j}\leq b_{\tau\left(  i\right)  }\right\}  $.
\par
Now, forget that we fixed $g$. We thus have shown that $g\in\left\{
j\in\left\{  1,2,\ldots,m\right\}  \ \mid\ a_{j}\leq b_{\tau\left(  i\right)
}\right\}  $ for each $g\in\left\{  j\in\left\{  1,2,\ldots,n\right\}
\ \mid\ b_{j}\leq b_{\tau\left(  i\right)  }\right\}  $. In other words,
$\left\{  j\in\left\{  1,2,\ldots,m\right\}  \ \mid\ b_{j}\leq b_{\tau\left(
i\right)  }\right\}  \subseteq\left\{  j\in\left\{  1,2,\ldots,n\right\}
\ \mid\ a_{j}\leq b_{\tau\left(  i\right)  }\right\}  $. Qed.}. Hence,%
\[
\left\vert \left\{  j\in\left\{  1,2,\ldots,m\right\}  \ \mid\ b_{j}\leq
b_{\tau\left(  i\right)  }\right\}  \right\vert \leq\left\vert \left\{
j\in\left\{  1,2,\ldots,n\right\}  \ \mid\ a_{j}\leq b_{\tau\left(  i\right)
}\right\}  \right\vert .
\]
Thus,%
\begin{align*}
&  \left\vert \left\{  j\in\left\{  1,2,\ldots,n\right\}  \ \mid\ a_{j}\leq
b_{\tau\left(  i\right)  }\right\}  \right\vert \\
&  \geq\left\vert \left\{  j\in\left\{  1,2,\ldots,m\right\}  \ \mid
\ b_{j}\leq b_{\tau\left(  i\right)  }\right\}  \right\vert \geq i.
\end{align*}
In other words, at least $i$ elements $j\in\left\{  1,2,\ldots,n\right\}  $
satisfy $a_{j}\leq b_{\tau\left(  i\right)  }$.
\end{verlong}

Now, we know that $b_{\tau\left(  i\right)  }$ is an element of $\mathbb{Z}$
such that at least $i$ elements $j\in\left\{  1,2,\ldots,n\right\}  $ satisfy
$a_{j}\leq b_{\tau\left(  i\right)  }$. In other words, $b_{\tau\left(
i\right)  }$ is an element $x$ of $\mathbb{Z}$ such that at least $i$ elements
$j\in\left\{  1,2,\ldots,n\right\}  $ satisfy $a_{j}\leq x$. In other words,%
\[
b_{\tau\left(  i\right)  }\in\left\{  x\in\mathbb{Z}\ \mid\ \text{at least
}i\text{ elements }j\in\left\{  1,2,\ldots,n\right\}  \text{ satisfy }%
a_{j}\leq x\right\}  .
\]
In light of (\ref{sol.sorting.nmu.X=}), this rewrites as $b_{\tau\left(
i\right)  }\in X$.

\begin{vershort}
But every $s\in X$ satisfies $\min X\leq s$ (since the minimum of a set is
$\leq$ to each element of this set). Applying this to $s=b_{\tau\left(
i\right)  }$, we obtain $\min X\leq b_{\tau\left(  i\right)  }$ (since
$b_{\tau\left(  i\right)  }\in X$). Thus, $b_{\tau\left(  i\right)  }\geq\min
X=a_{\sigma\left(  i\right)  }$ (by (\ref{sol.sorting.nmu.uniprop1'})). This
solves Exercise \ref{exe.sorting.nmu}. \qedhere

\end{vershort}

\begin{verlong}
But if $S$ is a subset of $\mathbb{Z}$ for which $\min S$ is well-defined,
then every $s\in S$ satisfies $\min S\leq s$ (since the minimum of a set is
$\leq$ to each element of this set). Applying this to $S=X$, we conclude that
every $s\in X$ satisfies $\min X\leq s$. Applying this to $s=b_{\tau\left(
i\right)  }$, we obtain $\min X\leq b_{\tau\left(  i\right)  }$ (since
$b_{\tau\left(  i\right)  }\in X$). Thus, $b_{\tau\left(  i\right)  }\geq\min
X=a_{\sigma\left(  i\right)  }$ (by (\ref{sol.sorting.nmu.uniprop1'})). In
other words, $a_{\sigma\left(  i\right)  }\leq b_{\tau\left(  i\right)  }$.
This solves Exercise \ref{exe.sorting.nmu}.
\end{verlong}
\end{proof}

In order to give a second solution to Exercise \ref{exe.sorting.nmu}, let us
prove a lemma:

\begin{lemma}
\label{lem.sorting.nmu.stronger}Let $n\in\mathbb{N}$ and $m\in\mathbb{N}$ be
such that $n\geq m$. Let $a_{1},a_{2},\ldots,a_{n}$ be $n$ integers. Let
$b_{1},b_{2},\ldots,b_{m}$ be $m$ integers. Assume that%
\begin{equation}
a_{i}\leq b_{i}\ \ \ \ \ \ \ \ \ \ \text{for every }i\in\left\{
1,2,\ldots,m\right\}  . \label{eq.lem.sorting.nmu.stronger.ass1}%
\end{equation}

Let $\sigma\in S_{n}$ and $\tau\in S_{m}$. Let $i\in\left\{  1,2,\ldots
,m\right\}  $. Assume that%
\begin{equation}
a_{\sigma\left(  i\right)  }\leq a_{\sigma\left(  v\right)  }%
\ \ \ \ \ \ \ \ \ \ \text{for every }v\in\left\{  i,i+1,\ldots,n\right\}  .
\label{eq.lem.sorting.nmu.stronger.ass2}%
\end{equation}
Furthermore, assume that%
\begin{equation}
b_{\tau\left(  u\right)  }\leq b_{\tau\left(  i\right)  }%
\ \ \ \ \ \ \ \ \ \ \text{for every }u\in\left\{  1,2,\ldots,i\right\}  .
\label{eq.lem.sorting.nmu.stronger.ass3}%
\end{equation}
Then, $a_{\sigma\left(  i\right)  }\leq b_{\tau\left(  i\right)  }$.
\end{lemma}

\begin{vershort}
\begin{proof}
[Proof of Lemma \ref{lem.sorting.nmu.stronger}.]We have $i\in\left\{
1,2,\ldots,m\right\}  \subseteq\left\{  1,2,\ldots,n\right\}  $ (since $m\leq
n$).

We have $\sigma\in S_{n}$. Thus, $\sigma$ is a permutation of $\left\{
1,2,\ldots,n\right\}  $, hence an injective map. Thus, the $n$ integers
$\sigma\left(  1\right)  ,\sigma\left(  2\right)  ,\ldots,\sigma\left(
n\right)  $ are distinct. In particular, the $n-i+1$ integers $\sigma\left(
i\right)  ,\sigma\left(  i+1\right)  ,\ldots,\sigma\left(  n\right)  $ are distinct.

Also, $\tau\in S_{m}$. Therefore, $\tau$ is a permutation of $\left\{
1,2,\ldots,m\right\}  $, hence an injective map. Therefore, the $m$ integers
$\tau\left(  1\right)  ,\tau\left(  2\right)  ,\ldots,\tau\left(  m\right)  $
are distinct. In particular, the $i$ integers $\tau\left(  1\right)
,\tau\left(  2\right)  ,\ldots,\tau\left(  i\right)  $ are distinct.

If $A$ and $B$ are two subsets of $\left\{  1,2,\ldots,n\right\}  $ satisfying
$\left\vert A\right\vert +\left\vert B\right\vert >n$, then%
\begin{equation}
A\cap B\neq\varnothing\label{pf.lem.sorting.nmu.stronger.short.AcutB}%
\end{equation}
\footnote{\textit{Proof of (\ref{pf.lem.sorting.nmu.stronger.short.AcutB}):}
Let $A$ and $B$ be two subsets of $\left\{  1,2,\ldots,n\right\}  $ satisfying
$\left\vert A\right\vert +\left\vert B\right\vert >n$. We must prove that
$A\cap B\neq\varnothing$.
\par
Indeed, assume the contrary. Thus, $A\cap B=\varnothing$. We know that $A$ and
$B$ are subsets of $\left\{  1,2,\ldots,n\right\}  $. Hence, $A\cup B$ is a
subset of $\left\{  1,2,\ldots,n\right\}  $. Thus, $\left\vert A\cup
B\right\vert \leq\left\vert \left\{  1,2,\ldots,n\right\}  \right\vert =n$.
But the sets $A$ and $B$ are disjoint (since $A\cap B=\varnothing$), and thus
we have $\left\vert A\cup B\right\vert =\left\vert A\right\vert +\left\vert
B\right\vert $ (since the size of the union of two disjoint sets is the sum of
their sizes). Thus, $\left\vert A\right\vert +\left\vert B\right\vert
=\left\vert A\cup B\right\vert \leq n$. This contradicts $\left\vert
A\right\vert +\left\vert B\right\vert >n$. This contradiction proves that our
assumption was wrong. Hence, $A\cap B\neq\varnothing$ is proven. This proves
(\ref{pf.lem.sorting.nmu.stronger.short.AcutB}).}.

Now, let $A=\left\{  \sigma\left(  i\right)  ,\sigma\left(  i+1\right)
,\ldots,\sigma\left(  n\right)  \right\}  $ and $B=\left\{  \tau\left(
1\right)  ,\tau\left(  2\right)  ,\ldots,\tau\left(  i\right)  \right\}  $.

Clearly,
\begin{equation}
A=\left\{  \sigma\left(  i\right)  ,\sigma\left(  i+1\right)  ,\ldots
,\sigma\left(  n\right)  \right\}  \subseteq\left\{  1,2,\ldots,n\right\}
\label{pf.lem.sorting.nmu.stronger.short.subset1}%
\end{equation}
(since $\sigma$ is a map $\left\{  1,2,\ldots,n\right\}  \rightarrow\left\{
1,2,\ldots,n\right\}  $) and $B=\left\{  \tau\left(  1\right)  ,\tau\left(
2\right)  ,\ldots,\tau\left(  i\right)  \right\}  \subseteq\left\{
1,2,\ldots,m\right\}  $ (since $\tau$ is a map $\left\{  1,2,\ldots,m\right\}
\rightarrow\left\{  1,2,\ldots,m\right\}  $). Thus,
\begin{equation}
B\subseteq\left\{  1,2,\ldots,m\right\}  \subseteq\left\{  1,2,\ldots
,n\right\}  . \label{pf.lem.sorting.nmu.stronger.short.subset2}%
\end{equation}

Now, from (\ref{pf.lem.sorting.nmu.stronger.short.subset1}) and
(\ref{pf.lem.sorting.nmu.stronger.short.subset2}), we see that $A$ and $B$ are
two subsets of $\left\{  1,2,\ldots,n\right\}  $. These two subsets satisfy%
\begin{align*}
&  \left\vert \underbrace{A}_{=\left\{  \sigma\left(  i\right)  ,\sigma\left(
i+1\right)  ,\ldots,\sigma\left(  n\right)  \right\}  }\right\vert +\left\vert
\underbrace{B}_{=\left\{  \tau\left(  1\right)  ,\tau\left(  2\right)
,\ldots,\tau\left(  i\right)  \right\}  }\right\vert \\
&  =\underbrace{\left\vert \left\{  \sigma\left(  i\right)  ,\sigma\left(
i+1\right)  ,\ldots,\sigma\left(  n\right)  \right\}  \right\vert
}_{\substack{=n-i+1\\\text{(since the }n-i+1\text{ integers}\\\sigma\left(
i\right)  ,\sigma\left(  i+1\right)  ,\ldots,\sigma\left(  n\right)  \text{
are distinct)}}}+\underbrace{\left\vert \left\{  \tau\left(  1\right)
,\tau\left(  2\right)  ,\ldots,\tau\left(  i\right)  \right\}  \right\vert
}_{\substack{=i\\\text{(since the }i\text{ integers}\\\tau\left(  1\right)
,\tau\left(  2\right)  ,\ldots,\tau\left(  i\right)  \text{ are distinct)}}}\\
&  =n-i+1+i=n+1>n.
\end{align*}
Hence, (\ref{pf.lem.sorting.nmu.stronger.short.AcutB}) shows that $A\cap
B\neq\varnothing$. In other words, there exists a $g\in A\cap B$. Consider
this $g$.

We have $g\in A\cap B\subseteq B=\left\{  \tau\left(  1\right)  ,\tau\left(
2\right)  ,\ldots,\tau\left(  i\right)  \right\}  \subseteq\left\{
1,2,\ldots,m\right\}  $. Hence, $b_{g}$ is well-defined. Also, $g\in\left\{
\tau\left(  1\right)  ,\tau\left(  2\right)  ,\ldots,\tau\left(  i\right)
\right\}  $. In other words, there exists a $u\in\left\{  1,2,\ldots
,i\right\}  $ satisfying $g=\tau\left(  u\right)  $. Consider this $u$. From
(\ref{eq.lem.sorting.nmu.stronger.ass3}), we obtain $b_{\tau\left(  u\right)
}\leq b_{\tau\left(  i\right)  }$. But $g=\tau\left(  u\right)  $ shows that
$b_{g}=b_{\tau\left(  u\right)  }\leq b_{\tau\left(  i\right)  }$.

On the other hand, $g\in A\cap B\subseteq A=\left\{  \sigma\left(  i\right)
,\sigma\left(  i+1\right)  ,\ldots,\sigma\left(  n\right)  \right\}
\subseteq\left\{  1,2,\ldots,n\right\}  $. Hence, $a_{g}$ is well-defined.
Also, $g\in\left\{  \sigma\left(  i\right)  ,\sigma\left(  i+1\right)
,\ldots,\sigma\left(  n\right)  \right\}  $. In other words, there exists a
$v\in\left\{  i,i+1,\ldots,n\right\}  $ such that $g=\sigma\left(  v\right)
$. Consider this $v$. From (\ref{eq.lem.sorting.nmu.stronger.ass2}), we obtain
$a_{\sigma\left(  i\right)  }\leq a_{\sigma\left(  v\right)  }=a_{g}$ (since
$\sigma\left(  v\right)  =g$).

But (\ref{eq.lem.sorting.nmu.stronger.ass1}) (applied to $g$ instead of $i$)
shows that $a_{g}\leq b_{g}$ (since $g\in\left\{  1,2,\ldots,m\right\}  $).
Hence, $a_{\sigma\left(  i\right)  }\leq a_{g}\leq b_{g}\leq b_{\tau\left(
i\right)  }$. This proves Lemma \ref{lem.sorting.nmu.stronger}.
\end{proof}
\end{vershort}

\begin{verlong}
\begin{proof}
[Proof of Lemma \ref{lem.sorting.nmu.stronger}.]We have $\sigma\in S_{n}$. In
other words, $\sigma$ belongs to the set $S_{n}$. In other words, $\sigma$
belongs to the set of all permutations of the set $\left\{  1,2,\ldots
,n\right\}  $ (since $S_{n}$ is the set of all permutations of the set
$\left\{  1,2,\ldots,n\right\}  $). In other words, $\sigma$ is a permutation
of the set $\left\{  1,2,\ldots,n\right\}  $. In other words, $\sigma$ is a
bijective map $\left\{  1,2,\ldots,n\right\}  \rightarrow\left\{
1,2,\ldots,n\right\}  $. The map $\sigma$ is thus bijective, hence both
surjective and injective.

We have $\tau\in S_{m}$. In other words, $\tau$ belongs to the set $S_{m}$. In
other words, $\tau$ belongs to the set of all permutations of the set
$\left\{  1,2,\ldots,m\right\}  $ (since $S_{m}$ is the set of all
permutations of the set $\left\{  1,2,\ldots,m\right\}  $). In other words,
$\tau$ is a permutation of the set $\left\{  1,2,\ldots,m\right\}  $. In other
words, $\tau$ is a bijective map $\left\{  1,2,\ldots,m\right\}
\rightarrow\left\{  1,2,\ldots,m\right\}  $. The map $\tau$ is thus bijective,
hence both surjective and injective.

We have $n\geq m$, thus $m\leq n$, hence $\left\{  1,2,\ldots,m\right\}
\subseteq\left\{  1,2,\ldots,n\right\}  $. Now, $i\in\left\{  1,2,\ldots
,m\right\}  \subseteq\left\{  1,2,\ldots,n\right\}  $; thus, $\sigma\left(
i\right)  $ is well-defined.

If $A$ and $B$ are two subsets of $\left\{  1,2,\ldots,n\right\}  $ satisfying
$\left\vert A\right\vert +\left\vert B\right\vert >n$, then%
\begin{equation}
A\cap B\neq\varnothing\label{pf.lem.sorting.nmu.stronger.AcutB}%
\end{equation}
\footnote{\textit{Proof of (\ref{pf.lem.sorting.nmu.stronger.AcutB}):} Let $A$
and $B$ be two subsets of $\left\{  1,2,\ldots,n\right\}  $ satisfying
$\left\vert A\right\vert +\left\vert B\right\vert >n$. We must prove that
$A\cap B\neq\varnothing$.
\par
Indeed, assume the contrary. Thus, $A\cap B=\varnothing$. We know that $A$ and
$B$ are subsets of $\left\{  1,2,\ldots,n\right\}  $. Hence, $A\cup B$ is a
subset of $\left\{  1,2,\ldots,n\right\}  $. In other words, $A\cup
B\subseteq\left\{  1,2,\ldots,n\right\}  $. Thus, $\left\vert A\cup
B\right\vert \leq\left\vert \left\{  1,2,\ldots,n\right\}  \right\vert =n$.
But the sets $A$ and $B$ are disjoint (since $A\cap B=\varnothing$), and thus
we have $\left\vert A\cup B\right\vert =\left\vert A\right\vert +\left\vert
B\right\vert $ (since the size of the union of two disjoint sets is the sum of
their sizes). Thus, $\left\vert A\right\vert +\left\vert B\right\vert
=\left\vert A\cup B\right\vert \leq n$. This contradicts $\left\vert
A\right\vert +\left\vert B\right\vert >n$. This contradiction proves that our
assumption was wrong. Hence, $A\cap B\neq\varnothing$ is proven. This proves
(\ref{pf.lem.sorting.nmu.stronger.AcutB}).}.

Now, let $A=\left\{  \sigma\left(  i\right)  ,\sigma\left(  i+1\right)
,\ldots,\sigma\left(  n\right)  \right\}  $ and $B=\left\{  \tau\left(
1\right)  ,\tau\left(  2\right)  ,\ldots,\tau\left(  i\right)  \right\}  $.

Clearly,
\begin{equation}
A=\left\{  \sigma\left(  i\right)  ,\sigma\left(  i+1\right)  ,\ldots
,\sigma\left(  n\right)  \right\}  \subseteq\left\{  1,2,\ldots,n\right\}
\label{pf.lem.sorting.nmu.stronger.subset1}%
\end{equation}
(since $\sigma$ is a map $\left\{  1,2,\ldots,n\right\}  \rightarrow\left\{
1,2,\ldots,n\right\}  $) and $B=\left\{  \tau\left(  1\right)  ,\tau\left(
2\right)  ,\ldots,\tau\left(  i\right)  \right\}  \subseteq\left\{
1,2,\ldots,m\right\}  $ (since $\tau$ is a map $\left\{  1,2,\ldots,m\right\}
\rightarrow\left\{  1,2,\ldots,m\right\}  $). Thus,
\begin{equation}
B\subseteq\left\{  1,2,\ldots,m\right\}  \subseteq\left\{  1,2,\ldots
,n\right\}  . \label{pf.lem.sorting.nmu.stronger.subset2}%
\end{equation}

The $n-i+1$ integers $\sigma\left(  i\right)  ,\sigma\left(  i+1\right)
,\ldots,\sigma\left(  n\right)  $ are pairwise
distinct\footnote{\textit{Proof.} Assume the contrary. Thus, the $n-i+1$
integers $\sigma\left(  i\right)  ,\sigma\left(  i+1\right)  ,\ldots
,\sigma\left(  n\right)  $ are \textbf{not} pairwise distinct. In other words,
two of these integers are equal. In other words, there exist two distinct
elements $u$ and $v$ of $\left\{  i,i+1,\ldots,n\right\}  $ such that
$\sigma\left(  u\right)  =\sigma\left(  v\right)  $. Consider these $u$ and
$v$.
\par
From $\sigma\left(  u\right)  =\sigma\left(  v\right)  $, we obtain $u=v$
(since the map $\sigma$ is injective). This contradicts the assumption that
$u$ and $v$ are distinct. This contradiction shows that our assumption was
wrong, qed.}. Hence,
\[
\left\vert \left\{  \sigma\left(  i\right)  ,\sigma\left(  i+1\right)
,\ldots,\sigma\left(  n\right)  \right\}  \right\vert =n-i+1.
\]

The $i$ integers $\tau\left(  1\right)  ,\tau\left(  2\right)  ,\ldots
,\tau\left(  i\right)  $ are pairwise distinct\footnote{\textit{Proof.} Assume
the contrary. Thus, the $i$ integers $\tau\left(  1\right)  ,\tau\left(
2\right)  ,\ldots,\tau\left(  i\right)  $ are \textbf{not} pairwise distinct.
In other words, two of these integers are equal. In other words, there exist
two distinct elements $u$ and $v$ of $\left\{  1,2,\ldots,i\right\}  $ such
that $\tau\left(  u\right)  =\tau\left(  v\right)  $. Consider these $u$ and
$v$.
\par
From $\tau\left(  u\right)  =\tau\left(  v\right)  $, we obtain $u=v$ (since
the map $\tau$ is injective). This contradicts the assumption that $u$ and $v$
are distinct. This contradiction shows that our assumption was wrong, qed.}.
Hence,
\[
\left\vert \left\{  \tau\left(  1\right)  ,\tau\left(  2\right)  ,\ldots
,\tau\left(  i\right)  \right\}  \right\vert =i.
\]

Now, from (\ref{pf.lem.sorting.nmu.stronger.subset1}) and
(\ref{pf.lem.sorting.nmu.stronger.subset2}), we see that $A$ and $B$ are two
subsets of $\left\{  1,2,\ldots,n\right\}  $. These two subsets satisfy%
\begin{align*}
&  \left\vert \underbrace{A}_{=\left\{  \sigma\left(  i\right)  ,\sigma\left(
i+1\right)  ,\ldots,\sigma\left(  n\right)  \right\}  }\right\vert +\left\vert
\underbrace{B}_{=\left\{  \tau\left(  1\right)  ,\tau\left(  2\right)
,\ldots,\tau\left(  i\right)  \right\}  }\right\vert \\
&  =\underbrace{\left\vert \left\{  \sigma\left(  i\right)  ,\sigma\left(
i+1\right)  ,\ldots,\sigma\left(  n\right)  \right\}  \right\vert }%
_{=n-i+1}+\underbrace{\left\vert \left\{  \tau\left(  1\right)  ,\tau\left(
2\right)  ,\ldots,\tau\left(  i\right)  \right\}  \right\vert }_{=i}\\
&  =n-i+1+i=n+1>n.
\end{align*}
Hence, (\ref{pf.lem.sorting.nmu.stronger.AcutB}) shows that $A\cap
B\neq\varnothing$. In other words, the set $A\cap B$ is nonempty. In other
words, the set $A\cap B$ has at least one element. In other words, there
exists a $g\in A\cap B$. Consider this $g$.

We have $g\in A\cap B\subseteq B=\left\{  \tau\left(  1\right)  ,\tau\left(
2\right)  ,\ldots,\tau\left(  i\right)  \right\}  \subseteq\left\{
1,2,\ldots,m\right\}  $. Hence, $b_{g}$ is well-defined. Also, $g\in\left\{
\tau\left(  1\right)  ,\tau\left(  2\right)  ,\ldots,\tau\left(  i\right)
\right\}  $. In other words, there exists a $u\in\left\{  1,2,\ldots
,i\right\}  $ satisfying $g=\tau\left(  u\right)  $. Consider this $u$. From
(\ref{eq.lem.sorting.nmu.stronger.ass3}), we obtain $b_{\tau\left(  u\right)
}\leq b_{\tau\left(  i\right)  }$. But $g=\tau\left(  u\right)  $ shows that
$b_{g}=b_{\tau\left(  u\right)  }\leq b_{\tau\left(  i\right)  }$.

On the other hand, $g\in A\cap B\subseteq A=\left\{  \sigma\left(  i\right)
,\sigma\left(  i+1\right)  ,\ldots,\sigma\left(  n\right)  \right\}
\subseteq\left\{  1,2,\ldots,n\right\}  $. Hence, $a_{g}$ is well-defined.
Also, $g\in\left\{  \sigma\left(  i\right)  ,\sigma\left(  i+1\right)
,\ldots,\sigma\left(  n\right)  \right\}  $. In other words, there exists a
$v\in\left\{  i,i+1,\ldots,n\right\}  $ such that $g=\sigma\left(  v\right)
$. Consider this $v$. From (\ref{eq.lem.sorting.nmu.stronger.ass2}), we obtain
$a_{\sigma\left(  i\right)  }\leq a_{\sigma\left(  v\right)  }=a_{g}$ (since
$\sigma\left(  v\right)  =g$).

But (\ref{eq.lem.sorting.nmu.stronger.ass1}) (applied to $g$ instead of $i$)
shows that $a_{g}\leq b_{g}$ (since $g\in\left\{  1,2,\ldots,m\right\}  $).
Hence, $a_{\sigma\left(  i\right)  }\leq a_{g}\leq b_{g}\leq b_{\tau\left(
i\right)  }$. This proves Lemma \ref{lem.sorting.nmu.stronger}.
\end{proof}
\end{verlong}

\begin{proof}
[Second solution to Exercise \ref{exe.sorting.nmu}.]Let $i\in\left\{
1,2,\ldots,m\right\}  $. We must prove that $a_{\sigma\left(  i\right)  }\leq
b_{\tau\left(  i\right)  }$.

We have $a_{\sigma\left(  1\right)  }\leq a_{\sigma\left(  2\right)  }%
\leq\cdots\leq a_{\sigma\left(  n\right)  }$. In other words, for every
$u\in\left\{  1,2,\ldots,n\right\}  $ and $v\in\left\{  1,2,\ldots,n\right\}
$ satisfying $u\leq v$, we have%
\begin{equation}
a_{\sigma\left(  u\right)  }\leq a_{\sigma\left(  v\right)  }.
\label{sol.sorting.nmu.sol2.1}%
\end{equation}
Also, we have $b_{\tau\left(  1\right)  }\leq b_{\tau\left(  2\right)  }%
\leq\cdots\leq b_{\tau\left(  m\right)  }$. In other words, for every
$u\in\left\{  1,2,\ldots,m\right\}  $ and $v\in\left\{  1,2,\ldots,m\right\}
$ satisfying $u\leq v$, we have%
\begin{equation}
b_{\tau\left(  u\right)  }\leq b_{\tau\left(  v\right)  }.
\label{sol.sorting.nmu.sol2.2}%
\end{equation}

Now,%
\[
a_{\sigma\left(  i\right)  }\leq a_{\sigma\left(  v\right)  }%
\ \ \ \ \ \ \ \ \ \ \text{for every }v\in\left\{  i,i+1,\ldots,n\right\}
\]
\footnote{\textit{Proof.} Let $v\in\left\{  i,i+1,\ldots,n\right\}  $. Then,
$v\geq i$, so that $i\leq v$. Also, $i\in\left\{  1,2,\ldots,m\right\}
\subseteq\left\{  1,2,\ldots,n\right\}  $ (since $m\leq n$ (since $n\geq m$))
and $v\in\left\{  i,i+1,\ldots,n\right\}  \subseteq\left\{  1,2,\ldots
,n\right\}  $ (since $i\geq1$ (since $i\in\left\{  1,2,\ldots,m\right\}  $)).
Hence, (\ref{sol.sorting.nmu.sol2.1}) (applied to $u=i$) yields $a_{\sigma
\left(  i\right)  }\leq a_{\sigma\left(  v\right)  }$. Qed.}. Furthermore,%
\[
b_{\tau\left(  u\right)  }\leq b_{\tau\left(  i\right)  }%
\ \ \ \ \ \ \ \ \ \ \text{for every }u\in\left\{  1,2,\ldots,i\right\}
\]
\footnote{\textit{Proof.} Let $u\in\left\{  1,2,\ldots,i\right\}  $. Then,
$u\leq i$. Also, $u\in\left\{  1,2,\ldots,i\right\}  \subseteq\left\{
1,2,\ldots,m\right\}  $ (since $i\leq m$ (since $i\in\left\{  1,2,\ldots
,m\right\}  $)) and $i\in\left\{  1,2,\ldots,m\right\}  $. Hence,
(\ref{sol.sorting.nmu.sol2.2}) (applied to $v=i$) yields $b_{\tau\left(
u\right)  }\leq b_{\tau\left(  i\right)  }$. Qed.}. Hence, Lemma
\ref{lem.sorting.nmu.stronger} shows that $a_{\sigma\left(  i\right)  }\leq
b_{\tau\left(  i\right)  }$. This solves Exercise \ref{exe.sorting.nmu} again.
\end{proof}

\subsection{Solution to Exercise \ref{exe.subsect.vandermonde.factoring}}

Definition \ref{def.vandermonde.factoring.h} shall be used throughout this section.

\begin{proof}
[Proof of Lemma \ref{lem.vandermonde.factoring.h.0}.]\textbf{(a)} Let $k$ be a
negative integer. Then,%
\begin{equation}
h_{k}\left(  x_{1},x_{2},\ldots,x_{n}\right)  =\sum_{\substack{\left(
a_{1},a_{2},\ldots,a_{n}\right)  \in\mathbb{N}^{n};\\a_{1}+a_{2}+\cdots
+a_{n}=k}}x_{1}^{a_{1}}x_{2}^{a_{2}}\cdots x_{n}^{a_{n}}
\label{pf.lem.vandermonde.factoring.h.0.a.1}%
\end{equation}
(by the definition of $h_{k}\left(  x_{1},x_{2},\ldots,x_{n}\right)  $).

\begin{vershort}
But $k<0$. Hence, there exists no $\left(  a_{1},a_{2},\ldots,a_{n}\right)
\in\mathbb{N}^{n}$ satisfying $a_{1}+a_{2}+\cdots+a_{n}=k$ (because every
$\left(  a_{1},a_{2},\ldots,a_{n}\right)  \in\mathbb{N}^{n}$ satisfies
$a_{1}+a_{2}+\cdots+a_{n}\geq0$, whereas $k<0$). Therefore, the sum on the
right hand side of (\ref{pf.lem.vandermonde.factoring.h.0.a.1}) is an empty
sum and therefore equals $0$. Therefore,
(\ref{pf.lem.vandermonde.factoring.h.0.a.1}) simplifies to $h_{k}\left(
x_{1},x_{2},\ldots,x_{n}\right)  =0$. This proves Lemma
\ref{lem.vandermonde.factoring.h.0} \textbf{(a)}.
\end{vershort}

\begin{verlong}
But there exists no $\left(  a_{1},a_{2},\ldots,a_{n}\right)  \in
\mathbb{N}^{n}$ satisfying $a_{1}+a_{2}+\cdots+a_{n}=k$%
\ \ \ \ \footnote{\textit{Proof.} Let $\left(  a_{1},a_{2},\ldots
,a_{n}\right)  \in\mathbb{N}^{n}$ be such that $a_{1}+a_{2}+\cdots+a_{n}=k$.
We shall derive a contradiction.
\par
We have $\left(  a_{1},a_{2},\ldots,a_{n}\right)  \in\mathbb{N}^{n}$. In other
words, $a_{1},a_{2},\ldots,a_{n}$ are $n$ elements of $\mathbb{N}$. In other
words, $a_{1},a_{2},\ldots,a_{n}$ are $n$ nonnegative integers. Hence, their
sum $a_{1}+a_{2}+\cdots+a_{n}$ also is a nonnegative integer. Thus,
$a_{1}+a_{2}+\cdots+a_{n}\geq0$. But $a_{1}+a_{2}+\cdots+a_{n}=k<0$ (since $k$
is negative). This contradicts $a_{1}+a_{2}+\cdots+a_{n}\geq0$.
\par
Now, forget that we fixed $\left(  a_{1},a_{2},\ldots,a_{n}\right)  $. We thus
have derived a contradiction for each $\left(  a_{1},a_{2},\ldots
,a_{n}\right)  \in\mathbb{N}^{n}$ satisfying $a_{1}+a_{2}+\cdots+a_{n}=k$.
Hence, there exists no $\left(  a_{1},a_{2},\ldots,a_{n}\right)  \in
\mathbb{N}^{n}$ satisfying $a_{1}+a_{2}+\cdots+a_{n}=k$. Qed.}. Hence, the sum
$\sum_{\substack{\left(  a_{1},a_{2},\ldots,a_{n}\right)  \in\mathbb{N}%
^{n};\\a_{1}+a_{2}+\cdots+a_{n}=k}}x_{1}^{a_{1}}x_{2}^{a_{2}}\cdots
x_{n}^{a_{n}}$ is an empty sum. Thus,
\[
\sum_{\substack{\left(  a_{1},a_{2},\ldots,a_{n}\right)  \in\mathbb{N}%
^{n};\\a_{1}+a_{2}+\cdots+a_{n}=k}}x_{1}^{a_{1}}x_{2}^{a_{2}}\cdots
x_{n}^{a_{n}}=\left(  \text{empty sum}\right)  =0.
\]
Hence, (\ref{pf.lem.vandermonde.factoring.h.0.a.1}) becomes%
\[
h_{k}\left(  x_{1},x_{2},\ldots,x_{n}\right)  =\sum_{\substack{\left(
a_{1},a_{2},\ldots,a_{n}\right)  \in\mathbb{N}^{n};\\a_{1}+a_{2}+\cdots
+a_{n}=k}}x_{1}^{a_{1}}x_{2}^{a_{2}}\cdots x_{n}^{a_{n}}=0.
\]
This proves Lemma \ref{lem.vandermonde.factoring.h.0} \textbf{(a)}.
\end{verlong}

\textbf{(b)} The definition of $h_{0}\left(  x_{1},x_{2},\ldots,x_{n}\right)
$ yields
\begin{equation}
h_{0}\left(  x_{1},x_{2},\ldots,x_{n}\right)  =\sum_{\substack{\left(
a_{1},a_{2},\ldots,a_{n}\right)  \in\mathbb{N}^{n};\\a_{1}+a_{2}+\cdots
+a_{n}=0}}x_{1}^{a_{1}}x_{2}^{a_{2}}\cdots x_{n}^{a_{n}}.
\label{pf.lem.vandermonde.factoring.h.0.b.1}%
\end{equation}

\begin{vershort}
But the only $\left(  a_{1},a_{2},\ldots,a_{n}\right)  \in\mathbb{N}^{n}$
satisfying $a_{1}+a_{2}+\cdots+a_{n}=0$ is $\left(  \underbrace{0,0,\ldots
,0}_{n\text{ zeroes}}\right)  $ (since a sum of nonnegative integers can only
be $0$ if all addends are $0$). Therefore, the sum on the right hand side of
(\ref{pf.lem.vandermonde.factoring.h.0.b.1}) has exactly one addend -- namely,
the addend for $\left(  a_{1},a_{2},\ldots,a_{n}\right)  =\left(
\underbrace{0,0,\ldots,0}_{n\text{ zeroes}}\right)  $. Thus,
(\ref{pf.lem.vandermonde.factoring.h.0.b.1}) simplifies to $h_{0}\left(
x_{1},x_{2},\ldots,x_{n}\right)  =\underbrace{x_{1}^{0}}_{=1}\underbrace{x_{2}%
^{0}}_{=1}\cdots\underbrace{x_{n}^{0}}_{=1}=1$. This proves Lemma
\ref{lem.vandermonde.factoring.h.0} \textbf{(b)}. \qedhere

\end{vershort}

\begin{verlong}
Define a subset $A$ of $\mathbb{N}^{n}$ by%
\begin{equation}
A=\left\{  \left(  a_{1},a_{2},\ldots,a_{n}\right)  \in\mathbb{N}^{n}%
\ \mid\ a_{1}+a_{2}+\cdots+a_{n}=0\right\}  .
\label{pf.lem.vandermonde.factoring.h.0.b.defA}%
\end{equation}
Thus,%
\begin{equation}
\sum_{\left(  a_{1},a_{2},\ldots,a_{n}\right)  \in A}=\sum_{\substack{\left(
a_{1},a_{2},\ldots,a_{n}\right)  \in\mathbb{N}^{n};\\a_{1}+a_{2}+\cdots
+a_{n}=0}} \label{pf.lem.vandermonde.factoring.h.0.b.sums1}%
\end{equation}
(an equality between summation signs). But $A=\left\{  \left(
\underbrace{0,0,\ldots,0}_{n\text{ zeroes}}\right)  \right\}  $%
\ \ \ \ \footnote{\textit{Proof.} The $n$-tuple $\left(
\underbrace{0,0,\ldots,0}_{n\text{ zeroes}}\right)  $ belongs to
$\mathbb{N}^{n}$ and satisfies $\underbrace{0+0+\cdots+0}_{n\text{ zeroes}}%
=0$. In other words, $\left(  \underbrace{0,0,\ldots,0}_{n\text{ zeroes}%
}\right)  $ is an element $\left(  a_{1},a_{2},\ldots,a_{n}\right)
\in\mathbb{N}^{n}$ satisfying $a_{1}+a_{2}+\cdots+a_{n}=0$. In other words,%
\[
\left(  \underbrace{0,0,\ldots,0}_{n\text{ zeroes}}\right)  \in\left\{
\left(  a_{1},a_{2},\ldots,a_{n}\right)  \in\mathbb{N}^{n}\ \mid\ a_{1}%
+a_{2}+\cdots+a_{n}=0\right\}  .
\]
In light of (\ref{pf.lem.vandermonde.factoring.h.0.b.defA}), this rewrites as
$\left(  \underbrace{0,0,\ldots,0}_{n\text{ zeroes}}\right)  \in A$. Thus,
$\left\{  \left(  \underbrace{0,0,\ldots,0}_{n\text{ zeroes}}\right)
\right\}  \subseteq A$.
\par
Now, let $\mathbf{a}\in A$. We shall show that $\mathbf{a}=\left(
\underbrace{0,0,\ldots,0}_{n\text{ zeroes}}\right)  $.
\par
Indeed, $\mathbf{a}\in A=\left\{  \left(  a_{1},a_{2},\ldots,a_{n}\right)
\in\mathbb{N}^{n}\ \mid\ a_{1}+a_{2}+\cdots+a_{n}=0\right\}  $. In other
words, $\mathbf{a}$ is an element $\left(  a_{1},a_{2},\ldots,a_{n}\right)
\in\mathbb{N}^{n}$ satisfying $a_{1}+a_{2}+\cdots+a_{n}=0$. In other words,
$\mathbf{a}=\left(  a_{1},a_{2},\ldots,a_{n}\right)  $ for some $\left(
a_{1},a_{2},\ldots,a_{n}\right)  \in\mathbb{N}^{n}$ satisfying $a_{1}%
+a_{2}+\cdots+a_{n}=0$. Consider this $\left(  a_{1},a_{2},\ldots
,a_{n}\right)  $.
\par
We have $\left(  a_{1},a_{2},\ldots,a_{n}\right)  \in\mathbb{N}^{n}$. In other
words, $a_{1},a_{2},\ldots,a_{n}$ are $n$ elements of $\mathbb{N}$.
\par
Fix $i\in\left\{  1,2,\ldots,n\right\}  $. We shall show that $a_{i}=0$.
\par
Splitting off the addend for $k=i$ from the sum $\sum_{k\in\left\{
1,2,\ldots,n\right\}  }a_{k}$, we obtain%
\[
\sum_{k\in\left\{  1,2,\ldots,n\right\}  }a_{k}=a_{i}+\sum_{\substack{k\in
\left\{  1,2,\ldots,n\right\}  ;\\k\neq i}}\underbrace{a_{k}}_{\substack{\geq
0\\\text{(since }a_{k}\in\mathbb{N}\\\text{(since }a_{1},a_{2},\ldots
,a_{n}\\\text{are }n\text{ elements of }\mathbb{N}\text{))}}}\geq
a_{i}+\underbrace{\sum_{\substack{k\in\left\{  1,2,\ldots,n\right\}  ;\\k\neq
i}}0}_{=0}=a_{i}.
\]
Thus, $a_{i}\leq\sum_{k\in\left\{  1,2,\ldots,n\right\}  }a_{k}=a_{1}%
+a_{2}+\cdots+a_{n}=0$. But $a_{i}\in\mathbb{N}$ (since $a_{1},a_{2}%
,\ldots,a_{n}$ are $n$ elements of $\mathbb{N}$), so that $a_{i}\geq0$.
Combining this with $a_{i}\leq0$, we obtain $a_{i}=0$.
\par
Now, forget that we fixed $i$. We thus have shown that $a_{i}=0$ for each
$i\in\left\{  1,2,\ldots,n\right\}  $. In other words, $\left(  a_{1}%
,a_{2},\ldots,a_{n}\right)  =\left(  \underbrace{0,0,\ldots,0}_{n\text{
zeroes}}\right)  $. Thus,
\[
\mathbf{a}=\left(  a_{1},a_{2},\ldots,a_{n}\right)  =\left(
\underbrace{0,0,\ldots,0}_{n\text{ zeroes}}\right)  \in\left\{  \left(
\underbrace{0,0,\ldots,0}_{n\text{ zeroes}}\right)  \right\}  .
\]
\par
Now, forget that we fixed $\mathbf{a}$. We thus have proven that
$\mathbf{a}\in\left\{  \left(  \underbrace{0,0,\ldots,0}_{n\text{ zeroes}%
}\right)  \right\}  $ for each $\mathbf{a}\in A$. In other words,
$A\subseteq\left\{  \left(  \underbrace{0,0,\ldots,0}_{n\text{ zeroes}%
}\right)  \right\}  $. Combining this with $\left\{  \left(
\underbrace{0,0,\ldots,0}_{n\text{ zeroes}}\right)  \right\}  \subseteq A$, we
obtain $A=\left\{  \left(  \underbrace{0,0,\ldots,0}_{n\text{ zeroes}}\right)
\right\}  $. Qed.}.

From $A=\left\{  \left(  \underbrace{0,0,\ldots,0}_{n\text{ zeroes}}\right)
\right\}  $, we obtain
\[
\sum_{\left(  a_{1},a_{2},\ldots,a_{n}\right)  \in A}=\sum_{\left(
a_{1},a_{2},\ldots,a_{n}\right)  \in\left\{  \left(  \underbrace{0,0,\ldots
,0}_{n\text{ zeroes}}\right)  \right\}  }%
\]
(an equality between summation signs). Comparing this with
(\ref{pf.lem.vandermonde.factoring.h.0.b.sums1}), we obtain%
\[
\sum_{\substack{\left(  a_{1},a_{2},\ldots,a_{n}\right)  \in\mathbb{N}%
^{n};\\a_{1}+a_{2}+\cdots+a_{n}=0}}=\sum_{\left(  a_{1},a_{2},\ldots
,a_{n}\right)  \in\left\{  \left(  \underbrace{0,0,\ldots,0}_{n\text{ zeroes}%
}\right)  \right\}  }%
\]
(an equality between summation signs). Now,
(\ref{pf.lem.vandermonde.factoring.h.0.b.1}) becomes%
\begin{align*}
h_{0}\left(  x_{1},x_{2},\ldots,x_{n}\right)   &  =\underbrace{\sum
_{\substack{\left(  a_{1},a_{2},\ldots,a_{n}\right)  \in\mathbb{N}^{n}%
;\\a_{1}+a_{2}+\cdots+a_{n}=0}}}_{=\sum_{\left(  a_{1},a_{2},\ldots
,a_{n}\right)  \in\left\{  \left(  \underbrace{0,0,\ldots,0}_{n\text{ zeroes}%
}\right)  \right\}  }}x_{1}^{a_{1}}x_{2}^{a_{2}}\cdots x_{n}^{a_{n}}\\
&  =\sum_{\left(  a_{1},a_{2},\ldots,a_{n}\right)  \in\left\{  \left(
\underbrace{0,0,\ldots,0}_{n\text{ zeroes}}\right)  \right\}  }x_{1}^{a_{1}%
}x_{2}^{a_{2}}\cdots x_{n}^{a_{n}}=x_{1}^{0}x_{2}^{0}\cdots x_{n}^{0}%
=\prod_{i=1}^{n}\underbrace{x_{i}^{0}}_{=1}\\
&  =\prod_{i=1}^{n}1=1.
\end{align*}
This proves Lemma \ref{lem.vandermonde.factoring.h.0} \textbf{(b)}.
\end{verlong}
\end{proof}

\begin{proof}
[Proof of Lemma \ref{lem.vandermonde.factoring.h.hq1}.]The definition of
$h_{q}\left(  x_{1},x_{2},\ldots,x_{k}\right)  $ yields%
\begin{align}
h_{q}\left(  x_{1},x_{2},\ldots,x_{k}\right)   &  =\sum_{\substack{\left(
a_{1},a_{2},\ldots,a_{k}\right)  \in\mathbb{N}^{k};\\a_{1}+a_{2}+\cdots
+a_{k}=q}}x_{1}^{a_{1}}x_{2}^{a_{2}}\cdots x_{k}^{a_{k}}\nonumber\\
&  =\sum_{\substack{\left(  s_{1},s_{2},\ldots,s_{k}\right)  \in\mathbb{N}%
^{k};\\s_{1}+s_{2}+\cdots+s_{k}=q}}x_{1}^{s_{1}}x_{2}^{s_{2}}\cdots
x_{k}^{s_{k}} \label{pf.lem.vandermonde.factoring.h.hq1.hq=}%
\end{align}
(here, we have renamed the summation index $\left(  a_{1},a_{2},\ldots
,a_{k}\right)  $ as $\left(  s_{1},s_{2},\ldots,s_{k}\right)  $).

Every $r\in\mathbb{Z}$ satisfies%
\begin{align}
h_{q-r}\left(  x_{1},x_{2},\ldots,x_{k-1}\right)   &  =\sum_{\substack{\left(
a_{1},a_{2},\ldots,a_{k-1}\right)  \in\mathbb{N}^{k-1};\\a_{1}+a_{2}%
+\cdots+a_{k-1}=q-r}}x_{1}^{a_{1}}x_{2}^{a_{2}}\cdots x_{k-1}^{a_{k-1}%
}\nonumber\\
&  \ \ \ \ \ \ \ \ \ \ \left(  \text{by the definition of }h_{q-r}\left(
x_{1},x_{2},\ldots,x_{k-1}\right)  \right) \nonumber\\
&  =\sum_{\substack{\left(  s_{1},s_{2},\ldots,s_{k-1}\right)  \in
\mathbb{N}^{k-1};\\s_{1}+s_{2}+\cdots+s_{k-1}=q-r}}x_{1}^{s_{1}}x_{2}^{s_{2}%
}\cdots x_{k-1}^{s_{k-1}} \label{pf.lem.vandermonde.factoring.h.hq1.hq-r=}%
\end{align}
(here, we have renamed the summation index $\left(  a_{1},a_{2},\ldots
,a_{k-1}\right)  $ as $\left(  s_{1},s_{2},\ldots,s_{k-1}\right)  $).
Moreover, every $r\in\mathbb{Z}$ satisfying $r>q$ satisfies%
\begin{equation}
h_{q-r}\left(  x_{1},x_{2},\ldots,x_{k-1}\right)  =0
\label{pf.lem.vandermonde.factoring.h.hq1.hq-r=0}%
\end{equation}
\footnote{\textit{Proof of (\ref{pf.lem.vandermonde.factoring.h.hq1.hq-r=0}):}
Let $r\in\mathbb{Z}$ be such that $r>q$. The integer $q-r$ is negative (since
$r>q$). Hence, Lemma \ref{lem.vandermonde.factoring.h.0} \textbf{(a)} (applied
to $k-1$ and $q-r$ instead of $n$ and $k$) yields $h_{q-r}\left(  x_{1}%
,x_{2},\ldots,x_{k-1}\right)  =0$. This proves
(\ref{pf.lem.vandermonde.factoring.h.hq1.hq-r=0}).}.

But Corollary \ref{cor.prodrule.prod-assMZ} (applied to $M=k$ and
$Z=\mathbb{N}$) shows that the map%
\begin{align*}
\mathbb{N}^{k}  &  \rightarrow\mathbb{N}^{k-1}\times\mathbb{N},\\
\left(  s_{1},s_{2},\ldots,s_{k}\right)   &  \mapsto\left(  \left(
s_{1},s_{2},\ldots,s_{k-1}\right)  ,s_{k}\right)
\end{align*}
is a bijection.

\begin{vershort}
Now, (\ref{pf.lem.vandermonde.factoring.h.hq1.hq=}) becomes%
\begin{align*}
h_{q}\left(  x_{1},x_{2},\ldots,x_{k}\right)   &  =\underbrace{\sum
_{\substack{\left(  s_{1},s_{2},\ldots,s_{k}\right)  \in\mathbb{N}^{k}%
;\\s_{1}+s_{2}+\cdots+s_{k}=q}}}_{\substack{=\sum_{\substack{\left(
s_{1},s_{2},\ldots,s_{k}\right)  \in\mathbb{N}^{k};\\s_{1}+s_{2}%
+\cdots+s_{k-1}=q-s_{k}}}\\\text{(because for any }\left(  s_{1},s_{2}%
,\ldots,s_{k}\right)  \in\mathbb{N}^{k}\text{,}\\\text{the condition }%
s_{1}+s_{2}+\cdots+s_{k}=q\\\text{is equivalent to }s_{1}+s_{2}+\cdots
+s_{k-1}=q-s_{k}\text{)}}}\underbrace{x_{1}^{s_{1}}x_{2}^{s_{2}}\cdots
x_{k}^{s_{k}}}_{=\left(  x_{1}^{s_{1}}x_{2}^{s_{2}}\cdots x_{k-1}^{s_{k-1}%
}\right)  x_{k}^{s_{k}}}\\
&  =\sum_{\substack{\left(  s_{1},s_{2},\ldots,s_{k}\right)  \in\mathbb{N}%
^{k};\\s_{1}+s_{2}+\cdots+s_{k-1}=q-s_{k}}}\left(  x_{1}^{s_{1}}x_{2}^{s_{2}%
}\cdots x_{k-1}^{s_{k-1}}\right)  x_{k}^{s_{k}}\\
&  =\underbrace{\sum_{\substack{\left(  \left(  s_{1},s_{2},\ldots
,s_{k-1}\right)  ,r\right)  \in\mathbb{N}^{k-1}\times\mathbb{N};\\s_{1}%
+s_{2}+\cdots+s_{k-1}=q-r}}}_{=\sum_{r\in\mathbb{N}}\sum_{\substack{\left(
s_{1},s_{2},\ldots,s_{k-1}\right)  \in\mathbb{N}^{k-1};\\s_{1}+s_{2}%
+\cdots+s_{k-1}=q-r}}}\left(  x_{1}^{s_{1}}x_{2}^{s_{2}}\cdots x_{k-1}%
^{s_{k-1}}\right)  x_{k}^{r}\\
&  \ \ \ \ \ \ \ \ \ \ \left(
\begin{array}
[c]{c}%
\text{here, we have substituted }\left(  \left(  s_{1},s_{2},\ldots
,s_{k-1}\right)  ,r\right)  \text{ for}\\
\left(  \left(  s_{1},s_{2},\ldots,s_{k-1}\right)  ,s_{k}\right)  \text{,
since the map}\\
\mathbb{N}^{k}\rightarrow\mathbb{N}^{k-1}\times\mathbb{N},\ \left(
s_{1},s_{2},\ldots,s_{k}\right)  \mapsto\left(  \left(  s_{1},s_{2}%
,\ldots,s_{k-1}\right)  ,s_{k}\right) \\
\text{is a bijection}%
\end{array}
\right)
\end{align*}%
\begin{align*}
&  =\sum_{r\in\mathbb{N}}\sum_{\substack{\left(  s_{1},s_{2},\ldots
,s_{k-1}\right)  \in\mathbb{N}^{k-1};\\s_{1}+s_{2}+\cdots+s_{k-1}=q-r}}\left(
x_{1}^{s_{1}}x_{2}^{s_{2}}\cdots x_{k-1}^{s_{k-1}}\right)  x_{k}^{r}\\
&  =\sum_{r\in\mathbb{N}}x_{k}^{r}\underbrace{\sum_{\substack{\left(
s_{1},s_{2},\ldots,s_{k-1}\right)  \in\mathbb{N}^{k-1};\\s_{1}+s_{2}%
+\cdots+s_{k-1}=q-r}}x_{1}^{s_{1}}x_{2}^{s_{2}}\cdots x_{k-1}^{s_{k-1}}%
}_{\substack{=h_{q-r}\left(  x_{1},x_{2},\ldots,x_{k-1}\right)  \\\text{(by
(\ref{pf.lem.vandermonde.factoring.h.hq1.hq-r=}))}}}\\
&  =\sum_{r\in\mathbb{N}}x_{k}^{r}h_{q-r}\left(  x_{1},x_{2},\ldots
,x_{k-1}\right) \\
&  =\underbrace{\sum_{\substack{r\in\mathbb{N};\\r\leq q}}}_{=\sum_{r=0}^{q}%
}x_{k}^{r}h_{q-r}\left(  x_{1},x_{2},\ldots,x_{k-1}\right)  +\sum
_{\substack{r\in\mathbb{N};\\r>q}}x_{k}^{r}\underbrace{h_{q-r}\left(
x_{1},x_{2},\ldots,x_{k-1}\right)  }_{\substack{=0\\\text{(by
(\ref{pf.lem.vandermonde.factoring.h.hq1.hq-r=0}))}}}\\
&  =\sum_{r=0}^{q}x_{k}^{r}h_{q-r}\left(  x_{1},x_{2},\ldots,x_{k-1}\right)
+\underbrace{\sum_{\substack{r\in\mathbb{N};\\r>q}}x_{k}^{r}0}_{=0}=\sum
_{r=0}^{q}x_{k}^{r}h_{q-r}\left(  x_{1},x_{2},\ldots,x_{k-1}\right)  .
\end{align*}
This proves Lemma \ref{lem.vandermonde.factoring.h.hq1}. \qedhere

\end{vershort}

\begin{verlong}
Now, (\ref{pf.lem.vandermonde.factoring.h.hq1.hq=}) becomes%
\begin{align*}
h_{q}\left(  x_{1},x_{2},\ldots,x_{k}\right)   &  =\underbrace{\sum
_{\substack{\left(  s_{1},s_{2},\ldots,s_{k}\right)  \in\mathbb{N}^{k}%
;\\s_{1}+s_{2}+\cdots+s_{k}=q}}}_{\substack{=\sum_{\substack{\left(
s_{1},s_{2},\ldots,s_{k}\right)  \in\mathbb{N}^{k};\\\left(  s_{1}%
+s_{2}+\cdots+s_{k-1}\right)  +s_{k}=q}}\\\text{(since }s_{1}+s_{2}%
+\cdots+s_{k}=\left(  s_{1}+s_{2}+\cdots+s_{k-1}\right)  +s_{k}\\\text{for
every }\left(  s_{1},s_{2},\ldots,s_{k}\right)  \in\mathbb{N}^{k}\text{)}%
}}\underbrace{x_{1}^{s_{1}}x_{2}^{s_{2}}\cdots x_{k}^{s_{k}}}_{=\left(
x_{1}^{s_{1}}x_{2}^{s_{2}}\cdots x_{k-1}^{s_{k-1}}\right)  x_{k}^{s_{k}}}\\
&  =\sum_{\substack{\left(  s_{1},s_{2},\ldots,s_{k}\right)  \in\mathbb{N}%
^{k};\\\left(  s_{1}+s_{2}+\cdots+s_{k-1}\right)  +s_{k}=q}}\left(
x_{1}^{s_{1}}x_{2}^{s_{2}}\cdots x_{k-1}^{s_{k-1}}\right)  x_{k}^{s_{k}}\\
&  =\underbrace{\sum_{\substack{\left(  \left(  s_{1},s_{2},\ldots
,s_{k-1}\right)  ,r\right)  \in\mathbb{N}^{k-1}\times\mathbb{N};\\\left(
s_{1}+s_{2}+\cdots+s_{k-1}\right)  +r=q}}}_{\substack{=\sum_{\left(
s_{1},s_{2},\ldots,s_{k-1}\right)  \in\mathbb{N}^{k-1}}\sum_{\substack{r\in
\mathbb{N};\\\left(  s_{1}+s_{2}+\cdots+s_{k-1}\right)  +r=q}}\\=\sum
_{r\in\mathbb{N}}\sum_{\substack{\left(  s_{1},s_{2},\ldots,s_{k-1}\right)
\in\mathbb{N}^{k-1};\\\left(  s_{1}+s_{2}+\cdots+s_{k-1}\right)  +r=q}%
}}}\left(  x_{1}^{s_{1}}x_{2}^{s_{2}}\cdots x_{k-1}^{s_{k-1}}\right)
x_{k}^{r}\\
&  \ \ \ \ \ \ \ \ \ \ \left(
\begin{array}
[c]{c}%
\text{here, we have substituted }\left(  \left(  s_{1},s_{2},\ldots
,s_{k-1}\right)  ,r\right)  \text{ for}\\
\left(  \left(  s_{1},s_{2},\ldots,s_{k-1}\right)  ,s_{k}\right)  \text{,
since the map}\\
\mathbb{N}^{k}\rightarrow\mathbb{N}^{k-1}\times\mathbb{N},\ \left(
s_{1},s_{2},\ldots,s_{k}\right)  \mapsto\left(  \left(  s_{1},s_{2}%
,\ldots,s_{k-1}\right)  ,s_{k}\right) \\
\text{is a bijection}%
\end{array}
\right) \\
&  =\sum_{r\in\mathbb{N}}\underbrace{\sum_{\substack{\left(  s_{1}%
,s_{2},\ldots,s_{k-1}\right)  \in\mathbb{N}^{k-1};\\\left(  s_{1}+s_{2}%
+\cdots+s_{k-1}\right)  +r=q}}}_{\substack{=\sum_{\substack{\left(
s_{1},s_{2},\ldots,s_{k-1}\right)  \in\mathbb{N}^{k-1};\\s_{1}+s_{2}%
+\cdots+s_{k-1}=q-r}}\\\text{(because for each }\left(  s_{1},s_{2}%
,\ldots,s_{k-1}\right)  \in\mathbb{N}^{k-1}\text{,}\\\text{the condition
}\left(  s_{1}+s_{2}+\cdots+s_{k-1}\right)  +r=q\\\text{is equivalent to
}s_{1}+s_{2}+\cdots+s_{k-1}=q-r\text{)}}}\left(  x_{1}^{s_{1}}x_{2}^{s_{2}%
}\cdots x_{k-1}^{s_{k-1}}\right)  x_{k}^{r}%
\end{align*}%
\begin{align*}
&  =\sum_{r\in\mathbb{N}}\sum_{\substack{\left(  s_{1},s_{2},\ldots
,s_{k-1}\right)  \in\mathbb{N}^{k-1};\\s_{1}+s_{2}+\cdots+s_{k-1}=q-r}}\left(
x_{1}^{s_{1}}x_{2}^{s_{2}}\cdots x_{k-1}^{s_{k-1}}\right)  x_{k}^{r}\\
&  =\sum_{r\in\mathbb{N}}x_{k}^{r}\underbrace{\sum_{\substack{\left(
s_{1},s_{2},\ldots,s_{k-1}\right)  \in\mathbb{N}^{k-1};\\s_{1}+s_{2}%
+\cdots+s_{k-1}=q-r}}x_{1}^{s_{1}}x_{2}^{s_{2}}\cdots x_{k-1}^{s_{k-1}}%
}_{\substack{=h_{q-r}\left(  x_{1},x_{2},\ldots,x_{k-1}\right)  \\\text{(by
(\ref{pf.lem.vandermonde.factoring.h.hq1.hq-r=}))}}}\\
&  =\sum_{r\in\mathbb{N}}x_{k}^{r}h_{q-r}\left(  x_{1},x_{2},\ldots
,x_{k-1}\right) \\
&  =\underbrace{\sum_{\substack{r\in\mathbb{N};\\r\leq q}}}_{=\sum_{r=0}^{q}%
}x_{k}^{r}h_{q-r}\left(  x_{1},x_{2},\ldots,x_{k-1}\right)  +\sum
_{\substack{r\in\mathbb{N};\\r>q}}x_{k}^{r}\underbrace{h_{q-r}\left(
x_{1},x_{2},\ldots,x_{k-1}\right)  }_{\substack{=0\\\text{(by
(\ref{pf.lem.vandermonde.factoring.h.hq1.hq-r=0}))}}}\\
&  \ \ \ \ \ \ \ \ \ \ \left(
\begin{array}
[c]{c}%
\text{since each }r\in\mathbb{N}\text{ satisfies either }r\leq q\text{ or
}r>q\\
\text{(but not both)}%
\end{array}
\right) \\
&  =\sum_{r=0}^{q}x_{k}^{r}h_{q-r}\left(  x_{1},x_{2},\ldots,x_{k-1}\right)
+\underbrace{\sum_{\substack{r\in\mathbb{N};\\r>q}}x_{k}^{r}0}_{=0}=\sum
_{r=0}^{q}x_{k}^{r}h_{q-r}\left(  x_{1},x_{2},\ldots,x_{k-1}\right)  .
\end{align*}
This proves Lemma \ref{lem.vandermonde.factoring.h.hq1}.
\end{verlong}
\end{proof}

\begin{vershort}
\begin{proof}
[Proof of Lemma \ref{lem.vandermonde.factoring.h.hq2}.]If $q<0$, then Lemma
\ref{lem.vandermonde.factoring.h.hq2} holds\footnote{\textit{Proof.} Assume
that $q<0$. Thus, both $q-1$ and $q$ are negative integers. Hence, Lemma
\ref{lem.vandermonde.factoring.h.0} \textbf{(a)} (applied to $q-1$ and $k$
instead of $k$ and $n$) yields $h_{q-1}\left(  x_{1},x_{2},\ldots
,x_{k}\right)  =0$. Also, Lemma \ref{lem.vandermonde.factoring.h.0}
\textbf{(a)} (applied to $q$ and $k-1$ instead of $k$ and $n$) yields
$h_{q}\left(  x_{1},x_{2},\ldots,x_{k-1}\right)  =0$.
\par
Now, Lemma \ref{lem.vandermonde.factoring.h.0} \textbf{(a)} (applied to $q$
and $k$ instead of $k$ and $n$) yields $h_{q}\left(  x_{1},x_{2},\ldots
,x_{k}\right)  =0$. Comparing this with%
\[
\underbrace{h_{q}\left(  x_{1},x_{2},\ldots,x_{k-1}\right)  }_{=0}%
+x_{k}\underbrace{h_{q-1}\left(  x_{1},x_{2},\ldots,x_{k}\right)  }%
_{=0}=0+x_{k}0=0,
\]
we obtain
\[
h_{q}\left(  x_{1},x_{2},\ldots,x_{k}\right)  =h_{q}\left(  x_{1},x_{2}%
,\ldots,x_{k-1}\right)  +x_{k}h_{q-1}\left(  x_{1},x_{2},\ldots,x_{k}\right)
.
\]
Thus, Lemma \ref{lem.vandermonde.factoring.h.hq2} holds. Qed.}. Hence, for the
rest of this proof, we WLOG assume that we don't have $q<0$. Thus, $q\geq0$.

Lemma \ref{lem.vandermonde.factoring.h.hq1} (applied to $q-1$ instead of $q$)
yields%
\begin{align}
h_{q-1}\left(  x_{1},x_{2},\ldots,x_{k}\right)   &  =\sum_{r=0}^{q-1}x_{k}%
^{r}h_{\left(  q-1\right)  -r}\left(  x_{1},x_{2},\ldots,x_{k-1}\right)
\nonumber\\
&  =\sum_{r=1}^{q}x_{k}^{r-1}\underbrace{h_{\left(  q-1\right)  -\left(
r-1\right)  }\left(  x_{1},x_{2},\ldots,x_{k-1}\right)  }_{\substack{=h_{q-r}%
\left(  x_{1},x_{2},\ldots,x_{k-1}\right)  \\\text{(since }\left(  q-1\right)
-\left(  r-1\right)  =q-r\text{)}}}\nonumber\\
&  \ \ \ \ \ \ \ \ \ \ \left(  \text{here, we have substituted }r-1\text{ for
}r\text{ in the sum}\right) \nonumber\\
&  =\sum_{r=1}^{q}x_{k}^{r-1}h_{q-r}\left(  x_{1},x_{2},\ldots,x_{k-1}\right)
. \label{pf.lem.vandermonde.factoring.h.hq2.short.1}%
\end{align}
But Lemma \ref{lem.vandermonde.factoring.h.hq1} yields%
\begin{align*}
h_{q}\left(  x_{1},x_{2},\ldots,x_{k}\right)   &  =\sum_{r=0}^{q}x_{k}%
^{r}h_{q-r}\left(  x_{1},x_{2},\ldots,x_{k-1}\right) \\
&  =\underbrace{x_{k}^{0}}_{=1}\underbrace{h_{q-0}\left(  x_{1},x_{2}%
,\ldots,x_{k-1}\right)  }_{=h_{q}\left(  x_{1},x_{2},\ldots,x_{k-1}\right)
}+\sum_{r=1}^{q}\underbrace{x_{k}^{r}}_{=x_{k}x_{k}^{r-1}}h_{q-r}\left(
x_{1},x_{2},\ldots,x_{k-1}\right) \\
&  \ \ \ \ \ \ \ \ \ \ \left(
\begin{array}
[c]{c}%
\text{here, we have split off the addend for }r=0\text{ from}\\
\text{the sum (since }q\geq0\text{)}%
\end{array}
\right) \\
&  =h_{q}\left(  x_{1},x_{2},\ldots,x_{k-1}\right)  +\underbrace{\sum
_{r=1}^{q}x_{k}x_{k}^{r-1}h_{q-r}\left(  x_{1},x_{2},\ldots,x_{k-1}\right)
}_{=x_{k}\sum_{r=1}^{q}x_{k}^{r-1}h_{q-r}\left(  x_{1},x_{2},\ldots
,x_{k-1}\right)  }\\
&  =h_{q}\left(  x_{1},x_{2},\ldots,x_{k-1}\right)  +x_{k}\underbrace{\sum
_{r=1}^{q}x_{k}^{r-1}h_{q-r}\left(  x_{1},x_{2},\ldots,x_{k-1}\right)
}_{\substack{=h_{q-1}\left(  x_{1},x_{2},\ldots,x_{k}\right)  \\\text{(by
(\ref{pf.lem.vandermonde.factoring.h.hq2.short.1}))}}}\\
&  =h_{q}\left(  x_{1},x_{2},\ldots,x_{k-1}\right)  +x_{k}h_{q-1}\left(
x_{1},x_{2},\ldots,x_{k}\right)  .
\end{align*}
This proves Lemma \ref{lem.vandermonde.factoring.h.hq2}.
\end{proof}
\end{vershort}

\begin{verlong}
\begin{proof}
[Proof of Lemma \ref{lem.vandermonde.factoring.h.hq2}.]If $q<0$, then Lemma
\ref{lem.vandermonde.factoring.h.hq2} holds\footnote{\textit{Proof.} Assume
that $q<0$. Thus, $q-1<q<0$ as well. In other words, $q-1$ is a negative
integer. Hence, Lemma \ref{lem.vandermonde.factoring.h.0} \textbf{(a)}
(applied to $q-1$ and $k$ instead of $k$ and $n$) yields $h_{q-1}\left(
x_{1},x_{2},\ldots,x_{k}\right)  =0$. Also, $q$ is a negative integer (since
$q<0$). Furthermore, $k-1\in\mathbb{N}$ (since $k$ is a positive integer).
Hence, Lemma \ref{lem.vandermonde.factoring.h.0} \textbf{(a)} (applied to $q$
and $k-1$ instead of $k$ and $n$) yields $h_{q}\left(  x_{1},x_{2}%
,\ldots,x_{k-1}\right)  =0$.
\par
Now, Lemma \ref{lem.vandermonde.factoring.h.0} \textbf{(a)} (applied to $q$
and $k$ instead of $k$ and $n$) yields $h_{q}\left(  x_{1},x_{2},\ldots
,x_{k}\right)  =0$. Comparing this with%
\[
\underbrace{h_{q}\left(  x_{1},x_{2},\ldots,x_{k-1}\right)  }_{=0}%
+x_{k}\underbrace{h_{q-1}\left(  x_{1},x_{2},\ldots,x_{k}\right)  }%
_{=0}=0+x_{k}0=0,
\]
we obtain
\[
h_{q}\left(  x_{1},x_{2},\ldots,x_{k}\right)  =h_{q}\left(  x_{1},x_{2}%
,\ldots,x_{k-1}\right)  +x_{k}h_{q-1}\left(  x_{1},x_{2},\ldots,x_{k}\right)
.
\]
Thus, Lemma \ref{lem.vandermonde.factoring.h.hq2} holds. Qed.}. Hence, for the
rest of this proof, we can WLOG assume that we don't have $q<0$. Assume this.

We have $q\geq0$ (since we don't have $q<0$). Thus, $q\in\mathbb{N}$.

Lemma \ref{lem.vandermonde.factoring.h.hq1} (applied to $q-1$ instead of $q$)
yields%
\begin{align}
h_{q-1}\left(  x_{1},x_{2},\ldots,x_{k}\right)   &  =\sum_{r=0}^{q-1}x_{k}%
^{r}h_{\left(  q-1\right)  -r}\left(  x_{1},x_{2},\ldots,x_{k-1}\right)
\nonumber\\
&  =\sum_{r=1}^{q}x_{k}^{r-1}\underbrace{h_{\left(  q-1\right)  -\left(
r-1\right)  }\left(  x_{1},x_{2},\ldots,x_{k-1}\right)  }_{\substack{=h_{q-r}%
\left(  x_{1},x_{2},\ldots,x_{k-1}\right)  \\\text{(since }\left(  q-1\right)
-\left(  r-1\right)  =q-r\text{)}}}\nonumber\\
&  \ \ \ \ \ \ \ \ \ \ \left(  \text{here, we have substituted }r-1\text{ for
}r\text{ in the sum}\right) \nonumber\\
&  =\sum_{r=1}^{q}x_{k}^{r-1}h_{q-r}\left(  x_{1},x_{2},\ldots,x_{k-1}\right)
. \label{pf.lem.vandermonde.factoring.h.hq2.1}%
\end{align}
But Lemma \ref{lem.vandermonde.factoring.h.hq1} yields%
\begin{align*}
h_{q}\left(  x_{1},x_{2},\ldots,x_{k}\right)   &  =\sum_{r=0}^{q}x_{k}%
^{r}h_{q-r}\left(  x_{1},x_{2},\ldots,x_{k-1}\right) \\
&  =\underbrace{x_{k}^{0}}_{=1}\underbrace{h_{q-0}\left(  x_{1},x_{2}%
,\ldots,x_{k-1}\right)  }_{\substack{=h_{q}\left(  x_{1},x_{2},\ldots
,x_{k-1}\right)  \\\text{(since }q-0=q\text{)}}}+\sum_{r=1}^{q}%
\underbrace{x_{k}^{r}}_{=x_{k}x_{k}^{r-1}}h_{q-r}\left(  x_{1},x_{2}%
,\ldots,x_{k-1}\right) \\
&  \ \ \ \ \ \ \ \ \ \ \left(
\begin{array}
[c]{c}%
\text{here, we have split off the addend for }r=0\text{ from the sum}\\
\text{(since }q\geq0\text{)}%
\end{array}
\right) \\
&  =h_{q}\left(  x_{1},x_{2},\ldots,x_{k-1}\right)  +\underbrace{\sum
_{r=1}^{q}x_{k}x_{k}^{r-1}h_{q-r}\left(  x_{1},x_{2},\ldots,x_{k-1}\right)
}_{=x_{k}\sum_{r=1}^{q}x_{k}^{r-1}h_{q-r}\left(  x_{1},x_{2},\ldots
,x_{k-1}\right)  }\\
&  =h_{q}\left(  x_{1},x_{2},\ldots,x_{k-1}\right)  +x_{k}\underbrace{\sum
_{r=1}^{q}x_{k}^{r-1}h_{q-r}\left(  x_{1},x_{2},\ldots,x_{k-1}\right)
}_{\substack{=h_{q-1}\left(  x_{1},x_{2},\ldots,x_{k}\right)  \\\text{(by
(\ref{pf.lem.vandermonde.factoring.h.hq2.1}))}}}\\
&  =h_{q}\left(  x_{1},x_{2},\ldots,x_{k-1}\right)  +x_{k}h_{q-1}\left(
x_{1},x_{2},\ldots,x_{k}\right)  .
\end{align*}
This proves Lemma \ref{lem.vandermonde.factoring.h.hq2}.
\end{proof}
\end{verlong}

\begin{proof}
[Proof of Lemma \ref{lem.vandermonde.factoring.h.sum-u}.]We first notice that
every positive integer $q$ satisfies%
\begin{equation}
h_{q}\left(  x_{1},x_{2},\ldots,x_{0}\right)  =0.
\label{pf.lem.vandermonde.factoring.h.sum-u.0}%
\end{equation}
(Here, as usual, $\left(  x_{1},x_{2},\ldots,x_{0}\right)  $ stands for the
$0$-tuple $\left(  {}\right)  $.)

\begin{vershort}
[\textit{Proof of (\ref{pf.lem.vandermonde.factoring.h.sum-u.0}):} Let $q$ be
a positive integer. The definition of $h_{q}\left(  x_{1},x_{2},\ldots
,x_{0}\right)  $ yields%
\begin{equation}
h_{q}\left(  x_{1},x_{2},\ldots,x_{0}\right)  =\sum_{\substack{\left(
a_{1},a_{2},\ldots,a_{0}\right)  \in\mathbb{N}^{0};\\a_{1}+a_{2}+\cdots
+a_{0}=q}}x_{1}^{a_{1}}x_{2}^{a_{2}}\cdots x_{0}^{a_{0}}.
\label{pf.lem.vandermonde.factoring.h.sum-u.0.pf.short.0}%
\end{equation}
But there is only one $0$-tuple $\left(  a_{1},a_{2},\ldots,a_{0}\right)
\in\mathbb{N}^{0}$, and this $0$-tuple satisfies $a_{1}+a_{2}+\cdots
+a_{0}=\left(  \text{empty sum}\right)  =0\neq q$ (since $q$ is positive).
Hence, there exists no $\left(  a_{1},a_{2},\ldots,a_{0}\right)  \in
\mathbb{N}^{0}$ satisfying $a_{1}+a_{2}+\cdots+a_{0}=q$. Thus, the sum on the
right hand side of (\ref{pf.lem.vandermonde.factoring.h.sum-u.0.pf.short.0})
is empty and equals $0$. Therefore,
(\ref{pf.lem.vandermonde.factoring.h.sum-u.0.pf.short.0}) simplifies to
$h_{q}\left(  x_{1},x_{2},\ldots,x_{0}\right)  =0$. This proves
(\ref{pf.lem.vandermonde.factoring.h.sum-u.0}).]
\end{vershort}

\begin{verlong}
[\textit{Proof of (\ref{pf.lem.vandermonde.factoring.h.sum-u.0}):} Let $q$ be
a positive integer. The definition of $h_{q}\left(  x_{1},x_{2},\ldots
,x_{0}\right)  $ yields%
\begin{equation}
h_{q}\left(  x_{1},x_{2},\ldots,x_{0}\right)  =\sum_{\substack{\left(
a_{1},a_{2},\ldots,a_{0}\right)  \in\mathbb{N}^{0};\\a_{1}+a_{2}+\cdots
+a_{0}=q}}x_{1}^{a_{1}}x_{2}^{a_{2}}\cdots x_{0}^{a_{0}}.
\label{pf.lem.vandermonde.factoring.h.sum-u.0.pf.0}%
\end{equation}
But there exists no $\left(  a_{1},a_{2},\ldots,a_{0}\right)  \in
\mathbb{N}^{0}$ satisfying $a_{1}+a_{2}+\cdots+a_{0}=q$%
\ \ \ \ \footnote{\textit{Proof.} We have $q\neq0$ (since $q$ is positive).
Every $\left(  a_{1},a_{2},\ldots,a_{0}\right)  \in\mathbb{N}^{0}$ satisfies
$a_{1}+a_{2}+\cdots+a_{0}=\left(  \text{empty sum}\right)  =0\neq q$. Hence,
no $\left(  a_{1},a_{2},\ldots,a_{0}\right)  \in\mathbb{N}^{0}$ satisfies
$a_{1}+a_{2}+\cdots+a_{0}=q$. In other words, there exists no $\left(
a_{1},a_{2},\ldots,a_{0}\right)  \in\mathbb{N}^{0}$ satisfying $a_{1}%
+a_{2}+\cdots+a_{0}=q$. Qed.}. Hence, the sum $\sum_{\substack{\left(
a_{1},a_{2},\ldots,a_{0}\right)  \in\mathbb{N}^{0};\\a_{1}+a_{2}+\cdots
+a_{0}=q}}x_{1}^{a_{1}}x_{2}^{a_{2}}\cdots x_{0}^{a_{0}}$ is empty. Thus,%
\[
\sum_{\substack{\left(  a_{1},a_{2},\ldots,a_{0}\right)  \in\mathbb{N}%
^{0};\\a_{1}+a_{2}+\cdots+a_{0}=q}}x_{1}^{a_{1}}x_{2}^{a_{2}}\cdots
x_{0}^{a_{0}}=\left(  \text{empty sum}\right)  =0.
\]
Thus, (\ref{pf.lem.vandermonde.factoring.h.sum-u.0.pf.0}) becomes%
\[
h_{q}\left(  x_{1},x_{2},\ldots,x_{0}\right)  =\sum_{\substack{\left(
a_{1},a_{2},\ldots,a_{0}\right)  \in\mathbb{N}^{0};\\a_{1}+a_{2}+\cdots
+a_{0}=q}}x_{1}^{a_{1}}x_{2}^{a_{2}}\cdots x_{0}^{a_{0}}=0.
\]
This proves (\ref{pf.lem.vandermonde.factoring.h.sum-u.0.pf.0}).]
\end{verlong}

We shall now show that%
\begin{equation}
\sum_{k=1}^{j}h_{j-k}\left(  x_{1},x_{2},\ldots,x_{k}\right)  \prod
_{p=1}^{k-1}\left(  u-x_{p}\right)  =u^{j-1}
\label{pf.lem.vandermonde.factoring.h.sum-u.main}%
\end{equation}
for each $j\in\left\{  1,2,\ldots,i\right\}  $.

[\textit{Proof of (\ref{pf.lem.vandermonde.factoring.h.sum-u.main}):} We shall
prove (\ref{pf.lem.vandermonde.factoring.h.sum-u.main}) by induction over $j$:

\begin{vershort}
\textit{Induction base:} We have%
\begin{align*}
\sum_{k=1}^{1}h_{1-k}\left(  x_{1},x_{2},\ldots,x_{k}\right)  \prod
_{p=1}^{k-1}\left(  u-x_{p}\right)   &  =\underbrace{h_{1-1}\left(
x_{1},x_{2},\ldots,x_{1}\right)  }_{\substack{=h_{0}\left(  x_{1},x_{2}%
,\ldots,x_{1}\right)  =1\\\text{(by Lemma \ref{lem.vandermonde.factoring.h.0}
\textbf{(b)}}\\\text{(applied to }1\text{ instead of }n\text{))}%
}}\underbrace{\prod_{p=1}^{1-1}\left(  u-x_{p}\right)  }_{=\left(  \text{empty
product}\right)  =1}\\
&  =1\cdot1=1=u^{1-1}%
\end{align*}
(since $u^{1-1}=u^{0}=1$). In other words,
(\ref{pf.lem.vandermonde.factoring.h.sum-u.main}) holds for $j=1$. This
completes the induction base.
\end{vershort}

\begin{verlong}
\textit{Induction base:} We have%
\begin{align*}
\sum_{k=1}^{1}h_{1-k}\left(  x_{1},x_{2},\ldots,x_{k}\right)  \prod
_{p=1}^{k-1}\left(  u-x_{p}\right)   &  =\underbrace{h_{1-1}\left(
x_{1},x_{2},\ldots,x_{1}\right)  }_{\substack{=h_{0}\left(  x_{1},x_{2}%
,\ldots,x_{1}\right)  \\\text{(since }1-1=0\text{)}}}\underbrace{\prod
_{p=1}^{1-1}\left(  u-x_{p}\right)  }_{\substack{=\left(  \text{empty
product}\right)  \\\text{(since }1-1<1\text{)}}}\\
&  =\underbrace{h_{0}\left(  x_{1},x_{2},\ldots,x_{1}\right)  }%
_{\substack{=1\\\text{(by Lemma \ref{lem.vandermonde.factoring.h.0}
\textbf{(b)}}\\\text{(applied to }1\text{ instead of }n\text{))}}%
}\cdot\underbrace{\left(  \text{empty product}\right)  }_{=1}\\
&  =1\cdot1=1.
\end{align*}
Comparing this with $u^{1-1}=u^{0}=1$, we obtain%
\[
\sum_{k=1}^{1}h_{1-k}\left(  x_{1},x_{2},\ldots,x_{k}\right)  \prod
_{p=1}^{k-1}\left(  u-x_{p}\right)  =u^{1-1}.
\]
In other words, (\ref{pf.lem.vandermonde.factoring.h.sum-u.main}) holds for
$j=1$. This completes the induction base.
\end{verlong}

\textit{Induction step:} Let $m\in\left\{  1,2,\ldots,i\right\}  $ be such
that $m>1$. Assume that (\ref{pf.lem.vandermonde.factoring.h.sum-u.main})
holds for $j=m-1$. We must prove that
(\ref{pf.lem.vandermonde.factoring.h.sum-u.main}) holds for $j=m$.

We have assumed that (\ref{pf.lem.vandermonde.factoring.h.sum-u.main}) holds
for $j=m-1$. In other words, we have%
\begin{equation}
\sum_{k=1}^{m-1}h_{\left(  m-1\right)  -k}\left(  x_{1},x_{2},\ldots
,x_{k}\right)  \prod_{p=1}^{k-1}\left(  u-x_{p}\right)  =u^{\left(
m-1\right)  -1}. \label{pf.lem.vandermonde.factoring.h.sum-u.main.pf.indhyp}%
\end{equation}

\begin{vershort}
Now, $m>1\geq0$. Moreover,
\begin{align}
&  \sum_{k=1}^{m}\underbrace{h_{m-k}\left(  x_{1},x_{2},\ldots,x_{k}\right)
}_{\substack{=h_{m-k}\left(  x_{1},x_{2},\ldots,x_{k-1}\right)  +x_{k}%
h_{m-k-1}\left(  x_{1},x_{2},\ldots,x_{k}\right)  \\\text{(by Lemma
\ref{lem.vandermonde.factoring.h.hq2} (applied to }q=m-k\text{))}}}\prod
_{p=1}^{k-1}\left(  u-x_{p}\right) \nonumber\\
&  =\sum_{k=1}^{m}\left(  h_{m-k}\left(  x_{1},x_{2},\ldots,x_{k-1}\right)
+x_{k}h_{m-k-1}\left(  x_{1},x_{2},\ldots,x_{k}\right)  \right)  \prod
_{p=1}^{k-1}\left(  u-x_{p}\right) \nonumber\\
&  =\sum_{k=1}^{m}h_{m-k}\left(  x_{1},x_{2},\ldots,x_{k-1}\right)
\prod_{p=1}^{k-1}\left(  u-x_{p}\right) \nonumber\\
&  \ \ \ \ \ \ \ \ \ \ +\sum_{k=1}^{m}x_{k}h_{m-k-1}\left(  x_{1},x_{2}%
,\ldots,x_{k}\right)  \prod_{p=1}^{k-1}\left(  u-x_{p}\right)  .
\label{pf.lem.vandermonde.factoring.h.sum-u.main.pf.short.1}%
\end{align}

\end{vershort}

\begin{verlong}
Now, $m>1\geq0$. Moreover,
\begin{align}
&  \sum_{k=1}^{m}\underbrace{h_{m-k}\left(  x_{1},x_{2},\ldots,x_{k}\right)
}_{\substack{=h_{m-k}\left(  x_{1},x_{2},\ldots,x_{k-1}\right)  +x_{k}%
h_{m-k-1}\left(  x_{1},x_{2},\ldots,x_{k}\right)  \\\text{(by Lemma
\ref{lem.vandermonde.factoring.h.hq2} (applied to }q=m-k\text{))}}}\prod
_{p=1}^{k-1}\left(  u-x_{p}\right) \nonumber\\
&  =\sum_{k=1}^{m}\underbrace{\left(  h_{m-k}\left(  x_{1},x_{2}%
,\ldots,x_{k-1}\right)  +x_{k}h_{m-k-1}\left(  x_{1},x_{2},\ldots
,x_{k}\right)  \right)  \prod_{p=1}^{k-1}\left(  u-x_{p}\right)  }%
_{=h_{m-k}\left(  x_{1},x_{2},\ldots,x_{k-1}\right)  \prod_{p=1}^{k-1}\left(
u-x_{p}\right)  +x_{k}h_{m-k-1}\left(  x_{1},x_{2},\ldots,x_{k}\right)
\prod_{p=1}^{k-1}\left(  u-x_{p}\right)  }\nonumber\\
&  =\sum_{k=1}^{m}\left(  h_{m-k}\left(  x_{1},x_{2},\ldots,x_{k-1}\right)
\prod_{p=1}^{k-1}\left(  u-x_{p}\right)  +x_{k}h_{m-k-1}\left(  x_{1}%
,x_{2},\ldots,x_{k}\right)  \prod_{p=1}^{k-1}\left(  u-x_{p}\right)  \right)
\nonumber\\
&  =\sum_{k=1}^{m}h_{m-k}\left(  x_{1},x_{2},\ldots,x_{k-1}\right)
\prod_{p=1}^{k-1}\left(  u-x_{p}\right) \nonumber\\
&  \ \ \ \ \ \ \ \ \ \ +\sum_{k=1}^{m}x_{k}h_{m-k-1}\left(  x_{1},x_{2}%
,\ldots,x_{k}\right)  \prod_{p=1}^{k-1}\left(  u-x_{p}\right)  .
\label{pf.lem.vandermonde.factoring.h.sum-u.main.pf.1}%
\end{align}

\end{verlong}

We have $h_{m-1}\left(  x_{1},x_{2},\ldots,x_{1-1}\right)  =0$%
\ \ \ \ \footnote{\textit{Proof.} We have $m>1$, and thus $m-1>0$. Hence,
$m-1$ is a positive integer. Thus,
(\ref{pf.lem.vandermonde.factoring.h.sum-u.0}) (applied to $q=m-1$) yields
$h_{m-1}\left(  x_{1},x_{2},\ldots,x_{0}\right)  =0$. Now, $1-1=0$ and thus
$h_{m-1}\left(  x_{1},x_{2},\ldots,x_{1-1}\right)  =h_{m-1}\left(  x_{1}%
,x_{2},\ldots,x_{0}\right)  =0$. Qed.}. But $m>1$. Hence, we can split off the
addend for $k=1$ from the sum $\sum_{k=1}^{m}h_{m-k}\left(  x_{1},x_{2}%
,\ldots,x_{k-1}\right)  \prod_{p=1}^{k-1}\left(  u-x_{p}\right)  $. We thus
obtain
\begin{align}
&  \sum_{k=1}^{m}h_{m-k}\left(  x_{1},x_{2},\ldots,x_{k-1}\right)  \prod
_{p=1}^{k-1}\left(  u-x_{p}\right) \nonumber\\
&  =\underbrace{h_{m-1}\left(  x_{1},x_{2},\ldots,x_{1-1}\right)  }_{=0}%
\prod_{p=1}^{1-1}\left(  u-x_{p}\right)  +\sum_{k=2}^{m}h_{m-k}\left(
x_{1},x_{2},\ldots,x_{k-1}\right)  \prod_{p=1}^{k-1}\left(  u-x_{p}\right)
\nonumber\\
&  =\underbrace{0\prod_{p=1}^{1-1}\left(  u-x_{p}\right)  }_{=0}+\sum
_{k=2}^{m}h_{m-k}\left(  x_{1},x_{2},\ldots,x_{k-1}\right)  \prod_{p=1}%
^{k-1}\left(  u-x_{p}\right) \nonumber\\
&  =\sum_{k=2}^{m}h_{m-k}\left(  x_{1},x_{2},\ldots,x_{k-1}\right)
\prod_{p=1}^{k-1}\left(  u-x_{p}\right) \nonumber\\
&  =\sum_{k=1}^{m-1}\underbrace{h_{m-\left(  k+1\right)  }\left(  x_{1}%
,x_{2},\ldots,x_{\left(  k+1\right)  -1}\right)  }_{\substack{=h_{m-k-1}%
\left(  x_{1},x_{2},\ldots,x_{k}\right)  \\\text{(since }m-\left(  k+1\right)
=m-k-1\\\text{and }\left(  k+1\right)  -1=k\text{)}}}\underbrace{\prod
_{p=1}^{\left(  k+1\right)  -1}}_{\substack{=\prod_{p=1}^{k}\\\text{(since
}\left(  k+1\right)  -1=k\text{)}}}\left(  u-x_{p}\right) \nonumber\\
&  \ \ \ \ \ \ \ \ \ \ \left(  \text{here, we have substituted }k+1\text{ for
}k\text{ in the sum}\right) \nonumber\\
&  =\sum_{k=1}^{m-1}h_{m-k-1}\left(  x_{1},x_{2},\ldots,x_{k}\right)
\underbrace{\prod_{p=1}^{k}\left(  u-x_{p}\right)  }_{\substack{=\left(
u-x_{k}\right)  \prod_{p=1}^{k-1}\left(  u-x_{p}\right)  \\\text{(here, we
have split off the factor for }p=k\text{ from}\\\text{the product (since }%
k\in\left\{  1,2,\ldots,k\right\}  \text{))}}}\nonumber\\
&  =\sum_{k=1}^{m-1}\underbrace{h_{m-k-1}\left(  x_{1},x_{2},\ldots
,x_{k}\right)  \cdot\left(  u-x_{k}\right)  }_{=\left(  u-x_{k}\right)
h_{m-k-1}\left(  x_{1},x_{2},\ldots,x_{k}\right)  }\prod_{p=1}^{k-1}\left(
u-x_{p}\right) \nonumber\\
&  =\sum_{k=1}^{m-1}\left(  u-x_{k}\right)  h_{m-k-1}\left(  x_{1}%
,x_{2},\ldots,x_{k}\right)  \prod_{p=1}^{k-1}\left(  u-x_{p}\right)  .
\label{pf.lem.vandermonde.factoring.h.sum-u.main.pf.2}%
\end{align}
On the other hand, $h_{m-m-1}\left(  x_{1},x_{2},\ldots,x_{m}\right)
=0$\ \ \ \ \footnote{\textit{Proof.} Clearly, $m-m-1=-1$ is a negative
integer. Hence, Lemma \ref{lem.vandermonde.factoring.h.0} \textbf{(a)}
(applied to $n=m$ and $k=m-m-1$) yields $h_{m-m-1}\left(  x_{1},x_{2}%
,\ldots,x_{m}\right)  =0$. Qed.}. But $m>1$. Hence, we can split off the
addend for $k=m$ from the sum $\sum_{k=1}^{m}x_{k}h_{m-k-1}\left(  x_{1}%
,x_{2},\ldots,x_{k}\right)  \prod_{p=1}^{k-1}\left(  u-x_{p}\right)  $. We
thus obtain%
\begin{align}
&  \sum_{k=1}^{m}x_{k}h_{m-k-1}\left(  x_{1},x_{2},\ldots,x_{k}\right)
\prod_{p=1}^{k-1}\left(  u-x_{p}\right) \nonumber\\
&  =x_{m}\underbrace{h_{m-m-1}\left(  x_{1},x_{2},\ldots,x_{m}\right)  }%
_{=0}\prod_{p=1}^{m-1}\left(  u-x_{p}\right)  +\sum_{k=1}^{m-1}x_{k}%
h_{m-k-1}\left(  x_{1},x_{2},\ldots,x_{k}\right)  \prod_{p=1}^{k-1}\left(
u-x_{p}\right) \nonumber\\
&  =\underbrace{x_{m}0\prod_{p=1}^{m-1}\left(  u-x_{p}\right)  }_{=0}%
+\sum_{k=1}^{m-1}x_{k}h_{m-k-1}\left(  x_{1},x_{2},\ldots,x_{k}\right)
\prod_{p=1}^{k-1}\left(  u-x_{p}\right) \nonumber\\
&  =\sum_{k=1}^{m-1}x_{k}h_{m-k-1}\left(  x_{1},x_{2},\ldots,x_{k}\right)
\prod_{p=1}^{k-1}\left(  u-x_{p}\right)  .
\label{pf.lem.vandermonde.factoring.h.sum-u.main.pf.3}%
\end{align}

\begin{vershort}
Now, (\ref{pf.lem.vandermonde.factoring.h.sum-u.main.pf.short.1}) becomes%
\begin{align*}
&  \sum_{k=1}^{m}h_{m-k}\left(  x_{1},x_{2},\ldots,x_{k}\right)  \prod
_{p=1}^{k-1}\left(  u-x_{p}\right) \\
&  =\underbrace{\sum_{k=1}^{m}h_{m-k}\left(  x_{1},x_{2},\ldots,x_{k-1}%
\right)  \prod_{p=1}^{k-1}\left(  u-x_{p}\right)  }_{\substack{=\sum
_{k=1}^{m-1}\left(  u-x_{k}\right)  h_{m-k-1}\left(  x_{1},x_{2},\ldots
,x_{k}\right)  \prod_{p=1}^{k-1}\left(  u-x_{p}\right)  \\\text{(by
(\ref{pf.lem.vandermonde.factoring.h.sum-u.main.pf.2}))}}}\\
&  \ \ \ \ \ \ \ \ \ \ +\underbrace{\sum_{k=1}^{m}x_{k}h_{m-k-1}\left(
x_{1},x_{2},\ldots,x_{k}\right)  \prod_{p=1}^{k-1}\left(  u-x_{p}\right)
}_{\substack{=\sum_{k=1}^{m-1}x_{k}h_{m-k-1}\left(  x_{1},x_{2},\ldots
,x_{k}\right)  \prod_{p=1}^{k-1}\left(  u-x_{p}\right)  \\\text{(by
(\ref{pf.lem.vandermonde.factoring.h.sum-u.main.pf.3}))}}}\\
&  =\sum_{k=1}^{m-1}\left(  u-x_{k}\right)  h_{m-k-1}\left(  x_{1}%
,x_{2},\ldots,x_{k}\right)  \prod_{p=1}^{k-1}\left(  u-x_{p}\right) \\
&  \ \ \ \ \ \ \ \ \ \ +\sum_{k=1}^{m-1}x_{k}h_{m-k-1}\left(  x_{1}%
,x_{2},\ldots,x_{k}\right)  \prod_{p=1}^{k-1}\left(  u-x_{p}\right) \\
&  =\sum_{k=1}^{m-1}\underbrace{\left(  \left(  u-x_{k}\right)  +x_{k}\right)
}_{=u}\underbrace{h_{m-k-1}\left(  x_{1},x_{2},\ldots,x_{k}\right)
}_{\substack{=h_{\left(  m-1\right)  -k}\left(  x_{1},x_{2},\ldots
,x_{k}\right)  \\\text{(since }m-k-1=\left(  m-1\right)  -k\text{)}}%
}\prod_{p=1}^{k-1}\left(  u-x_{p}\right) \\
&  =\sum_{k=1}^{m-1}uh_{\left(  m-1\right)  -k}\left(  x_{1},x_{2}%
,\ldots,x_{k}\right)  \prod_{p=1}^{k-1}\left(  u-x_{p}\right) \\
&  =u\underbrace{\sum_{k=1}^{m-1}h_{\left(  m-1\right)  -k}\left(  x_{1}%
,x_{2},\ldots,x_{k}\right)  \prod_{p=1}^{k-1}\left(  u-x_{p}\right)
}_{\substack{=u^{\left(  m-1\right)  -1}\\\text{(by
(\ref{pf.lem.vandermonde.factoring.h.sum-u.main.pf.indhyp}))}}}=uu^{\left(
m-1\right)  -1}=u^{m-1}.
\end{align*}
In other words, (\ref{pf.lem.vandermonde.factoring.h.sum-u.main}) holds for
$j=m$. This completes the induction step.
\end{vershort}

\begin{verlong}
Now, (\ref{pf.lem.vandermonde.factoring.h.sum-u.main.pf.1}) becomes%
\begin{align*}
&  \sum_{k=1}^{m}h_{m-k}\left(  x_{1},x_{2},\ldots,x_{k}\right)  \prod
_{p=1}^{k-1}\left(  u-x_{p}\right) \\
&  =\underbrace{\sum_{k=1}^{m}h_{m-k}\left(  x_{1},x_{2},\ldots,x_{k-1}%
\right)  \prod_{p=1}^{k-1}\left(  u-x_{p}\right)  }_{\substack{=\sum
_{k=1}^{m-1}\left(  u-x_{k}\right)  h_{m-k-1}\left(  x_{1},x_{2},\ldots
,x_{k}\right)  \prod_{p=1}^{k-1}\left(  u-x_{p}\right)  \\\text{(by
(\ref{pf.lem.vandermonde.factoring.h.sum-u.main.pf.2}))}}}\\
&  \ \ \ \ \ \ \ \ \ \ +\underbrace{\sum_{k=1}^{m}x_{k}h_{m-k-1}\left(
x_{1},x_{2},\ldots,x_{k}\right)  \prod_{p=1}^{k-1}\left(  u-x_{p}\right)
}_{\substack{=\sum_{k=1}^{m-1}x_{k}h_{m-k-1}\left(  x_{1},x_{2},\ldots
,x_{k}\right)  \prod_{p=1}^{k-1}\left(  u-x_{p}\right)  \\\text{(by
(\ref{pf.lem.vandermonde.factoring.h.sum-u.main.pf.3}))}}}\\
&  =\sum_{k=1}^{m-1}\left(  u-x_{k}\right)  h_{m-k-1}\left(  x_{1}%
,x_{2},\ldots,x_{k}\right)  \prod_{p=1}^{k-1}\left(  u-x_{p}\right) \\
&  \ \ \ \ \ \ \ \ \ \ +\sum_{k=1}^{m-1}x_{k}h_{m-k-1}\left(  x_{1}%
,x_{2},\ldots,x_{k}\right)  \prod_{p=1}^{k-1}\left(  u-x_{p}\right) \\
&  =\sum_{k=1}^{m-1}\underbrace{\left(  \left(  u-x_{k}\right)  h_{m-k-1}%
\left(  x_{1},x_{2},\ldots,x_{k}\right)  \prod_{p=1}^{k-1}\left(
u-x_{p}\right)  +x_{k}h_{m-k-1}\left(  x_{1},x_{2},\ldots,x_{k}\right)
\prod_{p=1}^{k-1}\left(  u-x_{p}\right)  \right)  }_{=\left(  \left(
u-x_{k}\right)  +x_{k}\right)  h_{m-k-1}\left(  x_{1},x_{2},\ldots
,x_{k}\right)  \prod_{p=1}^{k-1}\left(  u-x_{p}\right)  }\\
&  =\sum_{k=1}^{m-1}\underbrace{\left(  \left(  u-x_{k}\right)  +x_{k}\right)
}_{=u}\underbrace{h_{m-k-1}\left(  x_{1},x_{2},\ldots,x_{k}\right)
}_{\substack{=h_{\left(  m-1\right)  -k}\left(  x_{1},x_{2},\ldots
,x_{k}\right)  \\\text{(since }m-k-1=\left(  m-1\right)  -k\text{)}}%
}\prod_{p=1}^{k-1}\left(  u-x_{p}\right) \\
&  =\sum_{k=1}^{m-1}uh_{\left(  m-1\right)  -k}\left(  x_{1},x_{2}%
,\ldots,x_{k}\right)  \prod_{p=1}^{k-1}\left(  u-x_{p}\right) \\
&  =u\underbrace{\sum_{k=1}^{m-1}h_{\left(  m-1\right)  -k}\left(  x_{1}%
,x_{2},\ldots,x_{k}\right)  \prod_{p=1}^{k-1}\left(  u-x_{p}\right)
}_{\substack{=u^{\left(  m-1\right)  -1}\\\text{(by
(\ref{pf.lem.vandermonde.factoring.h.sum-u.main.pf.indhyp}))}}}=uu^{\left(
m-1\right)  -1}=u^{m-1}.
\end{align*}
In other words, (\ref{pf.lem.vandermonde.factoring.h.sum-u.main}) holds for
$j=m$. This completes the induction step.
\end{verlong}

Thus, we have proven (\ref{pf.lem.vandermonde.factoring.h.sum-u.main}) by induction.]

\begin{vershort}
Lemma \ref{lem.vandermonde.factoring.h.sum-u} now follows from
(\ref{pf.lem.vandermonde.factoring.h.sum-u.main}) (applied to $j=i$).
\qedhere

\end{vershort}

\begin{verlong}
Now, (\ref{pf.lem.vandermonde.factoring.h.sum-u.main}) (applied to $j=i$)
yields%
\[
\sum_{k=1}^{i}h_{i-k}\left(  x_{1},x_{2},\ldots,x_{k}\right)  \prod
_{p=1}^{k-1}\left(  u-x_{p}\right)  =u^{i-1}.
\]
This proves Lemma \ref{lem.vandermonde.factoring.h.sum-u}.
\end{verlong}
\end{proof}

\begin{proof}
[Proof of Lemma \ref{lem.vandermonde.factoring.U}.]For every $\left(
i,j\right)  \in\left\{  1,2,\ldots,n\right\}  ^{2}$, set
\begin{equation}
a_{i,j}=\prod_{p=1}^{j-1}\left(  x_{i}-x_{p}\right)  .
\label{pf.lem.vandermonde.factoring.U.1}%
\end{equation}
Recall that $U=\left(  \prod_{p=1}^{i-1}\left(  x_{j}-x_{p}\right)  \right)
_{1\leq i\leq n,\ 1\leq j\leq n}$. Hence, the definition of $U^{T}$ yields%
\[
U^{T}=\left(  \underbrace{\prod_{p=1}^{j-1}\left(  x_{i}-x_{p}\right)
}_{\substack{=a_{i,j}\\\text{(by (\ref{pf.lem.vandermonde.factoring.U.1}))}%
}}\right)  _{1\leq i\leq n,\ 1\leq j\leq n}=\left(  a_{i,j}\right)  _{1\leq
i\leq n,\ 1\leq j\leq n}.
\]

\begin{vershort}
We have $a_{i,j}=0$ for every $\left(  i,j\right)  \in\left\{  1,2,\ldots
,n\right\}  ^{2}$ satisfying $i<j$\ \ \ \ \footnote{\textit{Proof.} Let
$\left(  i,j\right)  \in\left\{  1,2,\ldots,n\right\}  ^{2}$ be such that
$i<j$. We want to show that $a_{i,j}=0$.
\par
We have $i<j$, and thus $i\leq j-1$ (since $i$ and $j$ are integers). Thus,
$i\in\left\{  1,2,\ldots,j-1\right\}  $. Thus, the product $\prod_{p=1}%
^{j-1}\left(  x_{i}-x_{p}\right)  $ has a factor for $p=i$. This factor is
$x_{i}-x_{i}=0$. Hence, at least one factor of the product $\prod_{p=1}%
^{j-1}\left(  x_{i}-x_{p}\right)  $ equals $0$ (namely, the factor for $p=i$).
Thus, the whole product $\prod_{p=1}^{j-1}\left(  x_{i}-x_{p}\right)  $ equals
$0$ (because if at least one factor of a product equals $0$, then the whole
product must equal $0$). In other words, $\prod_{p=1}^{j-1}\left(  x_{i}%
-x_{p}\right)  =0$.
\par
But the definition of $a_{i,j}$ yields $a_{i,j}=\prod_{p=1}^{j-1}\left(
x_{i}-x_{p}\right)  =0$. Qed.}. Hence, Exercise \ref{exe.ps4.3} (applied to
$U^{T}$ instead of $A$) yields
\begin{align*}
\det\left(  U^{T}\right)   &  =a_{1,1}a_{2,2}\cdots a_{n,n}=\prod_{i=1}%
^{n}\underbrace{a_{i,i}}_{\substack{=\prod_{p=1}^{i-1}\left(  x_{i}%
-x_{p}\right)  \\\text{(by the definition of }a_{i,i}\text{)}}%
}=\underbrace{\prod_{i=1}^{n}\prod_{p=1}^{i-1}}_{=\prod_{1\leq p<i\leq n}%
}\left(  x_{i}-x_{p}\right) \\
&  =\prod_{1\leq p<i\leq n}\left(  x_{i}-x_{p}\right)  =\prod_{1\leq j<i\leq
n}\left(  x_{i}-x_{j}\right)
\end{align*}
(here, we have renamed the index $\left(  p,i\right)  $ as $\left(
j,i\right)  $ in the product).
\end{vershort}

\begin{verlong}
We have $a_{i,j}=0$ for every $\left(  i,j\right)  \in\left\{  1,2,\ldots
,n\right\}  ^{2}$ satisfying $i<j$\ \ \ \ \footnote{\textit{Proof.} Let
$\left(  i,j\right)  \in\left\{  1,2,\ldots,n\right\}  ^{2}$ be such that
$i<j$. We want to show that $a_{i,j}=0$.
\par
We have $\left(  i,j\right)  \in\left\{  1,2,\ldots,n\right\}  ^{2}$. In other
words, $i\in\left\{  1,2,\ldots,n\right\}  $ and $j\in\left\{  1,2,\ldots
,n\right\}  $. From $i\in\left\{  1,2,\ldots,n\right\}  $, we obtain $i\geq1$.
\par
We have $i<j$, and thus $i\leq j-1$ (since $i$ and $j$ are integers). Thus,
$i\in\left\{  1,2,\ldots,j-1\right\}  $ (since $i\geq1$ and $i\leq j-1$).
Thus, the product $\prod_{p=1}^{j-1}\left(  x_{i}-x_{p}\right)  $ has a factor
for $p=i$. This factor is $x_{i}-x_{i}=0$. Hence, at least one factor of the
product $\prod_{p=1}^{j-1}\left(  x_{i}-x_{p}\right)  $ equals $0$ (namely,
the factor for $p=i$). Thus, the whole product $\prod_{p=1}^{j-1}\left(
x_{i}-x_{p}\right)  $ equals $0$ (because if at least one factor of a product
equals $0$, then the whole product must equal $0$). In other words,
$\prod_{p=1}^{j-1}\left(  x_{i}-x_{p}\right)  =0$.
\par
But the definition of $a_{i,j}$ yields $a_{i,j}=\prod_{p=1}^{j-1}\left(
x_{i}-x_{p}\right)  =0$. Qed.}. Hence, Exercise \ref{exe.ps4.3} (applied to
$U^{T}$ instead of $A$) yields
\begin{align*}
\det\left(  U^{T}\right)   &  =a_{1,1}a_{2,2}\cdots a_{n,n}=\prod_{i=1}%
^{n}\underbrace{a_{i,i}}_{\substack{=\prod_{p=1}^{i-1}\left(  x_{i}%
-x_{p}\right)  \\\text{(by the definition of }a_{i,i}\text{)}}%
}=\underbrace{\prod_{i=1}^{n}\prod_{p=1}^{i-1}}_{=\prod_{1\leq p<i\leq n}%
}\left(  x_{i}-x_{p}\right) \\
&  =\prod_{1\leq p<i\leq n}\left(  x_{i}-x_{p}\right)  =\prod_{1\leq j<i\leq
n}\left(  x_{i}-x_{j}\right)
\end{align*}
(here, we have renamed the index $\left(  p,i\right)  $ as $\left(
j,i\right)  $ in the product).
\end{verlong}

But Exercise \ref{exe.ps4.4} (applied to $A=U$) yields $\det\left(
U^{T}\right)  =\det U$. Hence,%
\[
\det U=\det\left(  U^{T}\right)  =\prod_{1\leq j<i\leq n}\left(  x_{i}%
-x_{j}\right)  .
\]
This proves Lemma \ref{lem.vandermonde.factoring.U}.
\end{proof}

\begin{proof}
[Proof of Lemma \ref{lem.vandermonde.factoring.L}.]For every $\left(
i,j\right)  \in\left\{  1,2,\ldots,n\right\}  ^{2}$, set
\begin{equation}
a_{i,j}=h_{i-j}\left(  x_{1},x_{2},\ldots,x_{j}\right)  .
\label{pf.lem.vandermonde.factoring.L.1}%
\end{equation}
Then,%
\[
L=\left(  \underbrace{h_{i-j}\left(  x_{1},x_{2},\ldots,x_{j}\right)
}_{\substack{=a_{i,j}\\\text{(by (\ref{pf.lem.vandermonde.factoring.L.1}))}%
}}\right)  _{1\leq i\leq n,\ 1\leq j\leq n}=\left(  a_{i,j}\right)  _{1\leq
i\leq n,\ 1\leq j\leq n}.
\]

We have $a_{i,j}=0$ for every $\left(  i,j\right)  \in\left\{  1,2,\ldots
,n\right\}  ^{2}$ satisfying $i<j$\ \ \ \ \footnote{\textit{Proof.} Let
$\left(  i,j\right)  \in\left\{  1,2,\ldots,n\right\}  ^{2}$ be such that
$i<j$. We want to show that $a_{i,j}=0$.
\par
We have $i-j<0$ (since $i<j$). Thus, $i-j$ is a negative integer. Hence, Lemma
\ref{lem.vandermonde.factoring.h.0} \textbf{(a)} (applied to $j$ and $i-j$
instead of $n$ and $k$) yields $h_{i-j}\left(  x_{1},x_{2},\ldots
,x_{j}\right)  =0$. But the definition of $a_{i,j}$ yields $a_{i,j}%
=h_{i-j}\left(  x_{1},x_{2},\ldots,x_{j}\right)  =0$. Qed.}. Hence, Exercise
\ref{exe.ps4.3} (applied to $L$ instead of $A$) yields
\begin{align*}
\det L  &  =a_{1,1}a_{2,2}\cdots a_{n,n}=\prod_{i=1}^{n}\underbrace{a_{i,i}%
}_{\substack{=h_{i-i}\left(  x_{1},x_{2},\ldots,x_{i}\right)  \\\text{(by the
definition of }a_{i,i}\text{)}}}=\prod_{i=1}^{n}\underbrace{h_{i-i}\left(
x_{1},x_{2},\ldots,x_{i}\right)  }_{\substack{=h_{0}\left(  x_{1},x_{2}%
,\ldots,x_{i}\right)  \\\text{(since }i-i=0\text{)}}}\\
&  =\prod_{i=1}^{n}\underbrace{h_{0}\left(  x_{1},x_{2},\ldots,x_{i}\right)
}_{\substack{=1\\\text{(by Lemma \ref{lem.vandermonde.factoring.h.0}
\textbf{(b)}}\\\text{(applied to }i\text{ instead of }n\text{))}}}=\prod
_{i=1}^{n}1=1.
\end{align*}
This proves Lemma \ref{lem.vandermonde.factoring.L}.
\end{proof}

\begin{proof}
[Proof of Lemma \ref{lem.vandermonde.factoring.A=LU}.]If $i\in\mathbb{Z}$ and
$k\in\left\{  1,2,\ldots,n\right\}  $ are such that $k>i$, then%
\begin{equation}
h_{i-k}\left(  x_{1},x_{2},\ldots,x_{k}\right)  =0
\label{pf.lem.vandermonde.factoring.A=LU.0}%
\end{equation}
\footnote{\textit{Proof of (\ref{pf.lem.vandermonde.factoring.A=LU.0}):} Let
$i\in\mathbb{Z}$ and $k\in\left\{  1,2,\ldots,n\right\}  $ be such that $k>i$.
Then, $i-k<0$ (since $k>i$). Thus, $i-k$ is a negative integer. Hence, Lemma
\ref{lem.vandermonde.factoring.h.0} \textbf{(a)} (applied to $k$ and $i-k$
instead of $n$ and $k$) yields $h_{i-k}\left(  x_{1},x_{2},\ldots
,x_{k}\right)  =0$. This proves (\ref{pf.lem.vandermonde.factoring.A=LU.0}).}.

Now, for every $i\in\left\{  1,2,\ldots,n\right\}  $ and $u\in\mathbb{K}$, we
have%
\begin{align}
&  \sum_{k=1}^{n}h_{i-k}\left(  x_{1},x_{2},\ldots,x_{k}\right)  \prod
_{p=1}^{k-1}\left(  u-x_{p}\right) \nonumber\\
&  =\underbrace{\sum_{k=1}^{i}h_{i-k}\left(  x_{1},x_{2},\ldots,x_{k}\right)
\prod_{p=1}^{k-1}\left(  u-x_{p}\right)  }_{\substack{=u^{i-1}\\\text{(by
Lemma \ref{lem.vandermonde.factoring.h.sum-u})}}}+\sum_{k=i+1}^{n}%
\underbrace{h_{i-k}\left(  x_{1},x_{2},\ldots,x_{k}\right)  }%
_{\substack{=0\\\text{(by (\ref{pf.lem.vandermonde.factoring.A=LU.0}%
)}\\\text{(since }k\geq i+1>i\text{))}}}\prod_{p=1}^{k-1}\left(
u-x_{p}\right) \nonumber\\
&  \ \ \ \ \ \ \ \ \ \ \left(  \text{since }1\leq i\leq n\text{ (since }%
i\in\left\{  1,2,\ldots,n\right\}  \text{)}\right) \nonumber\\
&  =u^{i-1}+\underbrace{\sum_{k=i+1}^{n}0\prod_{p=1}^{k-1}\left(
u-x_{p}\right)  }_{=0}=u^{i-1}. \label{pf.lem.vandermonde.factoring.A=LU.2}%
\end{align}

But recall that $L=\left(  h_{i-j}\left(  x_{1},x_{2},\ldots,x_{j}\right)
\right)  _{1\leq i\leq n,\ 1\leq j\leq n}$ and \newline$U=\left(  \prod
_{p=1}^{i-1}\left(  x_{j}-x_{p}\right)  \right)  _{1\leq i\leq n,\ 1\leq j\leq
n}$. Hence, the definition of the product $LU$ yields%
\[
LU=\left(  \underbrace{\sum_{k=1}^{n}h_{i-k}\left(  x_{1},x_{2},\ldots
,x_{k}\right)  \prod_{p=1}^{k-1}\left(  x_{j}-x_{p}\right)  }%
_{\substack{=x_{j}^{i-1}\\\text{(by (\ref{pf.lem.vandermonde.factoring.A=LU.2}%
) (applied to }u=x_{j}\text{))}}}\right)  _{1\leq i\leq n,\ 1\leq j\leq
n}=\left(  x_{j}^{i-1}\right)  _{1\leq i\leq n,\ 1\leq j\leq n}.
\]
This proves Lemma \ref{lem.vandermonde.factoring.A=LU}.
\end{proof}

\begin{proof}
[Solution to Exercise \ref{exe.subsect.vandermonde.factoring}.]We have proven
Lemma \ref{lem.vandermonde.factoring.h.0}, Lemma
\ref{lem.vandermonde.factoring.h.hq1}, Lemma
\ref{lem.vandermonde.factoring.h.hq2}, Lemma
\ref{lem.vandermonde.factoring.h.sum-u}, Lemma
\ref{lem.vandermonde.factoring.U}, Lemma \ref{lem.vandermonde.factoring.L} and
Lemma \ref{lem.vandermonde.factoring.A=LU}. Thus, Exercise
\ref{exe.subsect.vandermonde.factoring} is solved.
\end{proof}

\subsection{Solution to Exercise \ref{exe.vander-det.s1}}

\begin{proof}
[First solution to Exercise \ref{exe.vander-det.s1}.]Our solution will imitate
our First proof of Theorem \ref{thm.vander-det} \textbf{(a)} (but it will
involve some additional complications).

For every $u\in\left\{  0,1,\ldots,n\right\}  $ and $\left(  i,j\right)
\in\left\{  1,2,\ldots,u\right\}  ^{2}$, define $a_{i,j,u}\in\mathbb{K}$ by%
\[
a_{i,j,u}=%
\begin{cases}
x_{i}^{u-j}, & \text{if }j>1;\\
x_{i}^{u}, & \text{if }j=1
\end{cases}
.
\]

For every $u\in\left\{  1,2,\ldots,n\right\}  $, let $A_{u}$ be the $u\times
u$-matrix $\left(  a_{i,j,u}\right)  _{1\leq i\leq u,\ 1\leq j\leq u}$. Then,
our goal is to prove that $\det\left(  A_{n}\right)  =\left(  x_{1}%
+x_{2}+\cdots+x_{n}\right)  \prod_{1\leq i<j\leq n}\left(  x_{i}-x_{j}\right)
$ (because $A_{n}$ is precisely the matrix $\left(
\begin{cases}
x_{i}^{n-j}, & \text{if }j>1;\\
x_{i}^{n}, & \text{if }j=1
\end{cases}
\right)  _{1\leq i\leq n,\ 1\leq j\leq n}$ which appears in the statement of
the exercise).

Now, let us show that%
\begin{equation}
\det\left(  A_{u}\right)  =\left(  x_{1}+x_{2}+\cdots+x_{u}\right)
\prod_{1\leq i<j\leq u}\left(  x_{i}-x_{j}\right)
\label{sol.vander-det.s1.goal}%
\end{equation}
for every $u\in\left\{  1,2,\ldots,n\right\}  $.

[\textit{Proof of (\ref{sol.vander-det.s1.goal}):} We will prove
(\ref{sol.vander-det.s1.goal}) by induction over $u$:

\textit{Induction base:} The definition of $A_{1}$ yields $A_{1}=\left(
\begin{array}
[c]{c}%
x_{1}^{1}%
\end{array}
\right)  =\left(
\begin{array}
[c]{c}%
x_{1}%
\end{array}
\right)  $ and thus $\det\left(  A_{1}\right)  =x_{1}$. Compared with%
\[
\underbrace{\left(  x_{1}+x_{2}+\cdots+x_{1}\right)  }_{=x_{1}}%
\underbrace{\prod_{1\leq i<j\leq1}\left(  x_{i}-x_{j}\right)  }_{=\left(
\text{empty product}\right)  =1}=x_{1},
\]
this yields $\det\left(  A_{1}\right)  =\left(  x_{1}+x_{2}+\cdots
+x_{1}\right)  \prod_{1\leq i<j\leq1}\left(  x_{i}-x_{j}\right)  $. In other
words, (\ref{sol.vander-det.s1.goal}) holds for $u=1$. The induction base is
thus complete.

\textit{Induction step:} Let $U\in\left\{  2,3,\ldots,n\right\}  $. Assume
that (\ref{sol.vander-det.s1.goal}) holds for $u=U-1$. We need to prove that
(\ref{sol.vander-det.s1.goal}) holds for $u=U$.

Recall that $A_{U}=\left(  a_{i,j,U}\right)  _{1\leq i\leq U,\ 1\leq j\leq U}$
(by the definition of $A_{U}$). The matrix $A_{U}$ looks as follows:%
\[
A_{U}=\left(
\begin{array}
[c]{cccccc}%
x_{1}^{U} & x_{1}^{U-2} & x_{1}^{U-3} & \cdots & x_{1} & 1\\
x_{2}^{U} & x_{2}^{U-2} & x_{2}^{U-3} & \cdots & x_{2} & 1\\
x_{3}^{U} & x_{3}^{U-2} & x_{3}^{U-3} & \cdots & x_{3} & 1\\
x_{4}^{U} & x_{4}^{U-2} & x_{4}^{U-3} & \cdots & x_{4} & 1\\
\vdots & \vdots & \vdots & \ddots & \vdots & \vdots\\
x_{U}^{U} & x_{U}^{U-2} & x_{U}^{U-3} & \cdots & x_{U} & 1
\end{array}
\right)  .
\]

For every $\left(  i,j\right)  \in\left\{  1,2,\ldots,U\right\}  ^{2}$, define
$b_{i,j}\in\mathbb{K}$ by%
\[
b_{i,j}=%
\begin{cases}
x_{i}^{U}-x_{U}^{2}x_{i}^{U-2}, & \text{if }j=1;\\
x_{i}^{U-j}-x_{U}x_{i}^{U-j-1}, & \text{if }1<j<U;\\
1, & \text{if }j=U
\end{cases}
.
\]
Let $B$ be the $U\times U$-matrix $\left(  b_{i,j}\right)  _{1\leq i\leq
U,\ 1\leq j\leq U}$. Here is how $B$ looks like:%
\[
B=\left(
\begin{array}
[c]{cccccc}%
x_{1}^{U}-x_{U}^{2}x_{1}^{U-2} & x_{1}^{U-2}-x_{U}x_{1}^{U-3} & x_{1}%
^{U-3}-x_{U}x_{1}^{U-4} & \cdots & x_{1}-x_{U} & 1\\
x_{2}^{U}-x_{U}^{2}x_{2}^{U-2} & x_{2}^{U-2}-x_{U}x_{2}^{U-3} & x_{2}%
^{U-3}-x_{U}x_{2}^{U-4} & \cdots & x_{2}-x_{U} & 1\\
x_{3}^{U}-x_{U}^{2}x_{3}^{U-2} & x_{3}^{U-2}-x_{U}x_{3}^{U-3} & x_{3}%
^{U-3}-x_{U}x_{3}^{U-4} & \cdots & x_{3}-x_{U} & 1\\
x_{4}^{U}-x_{U}^{2}x_{4}^{U-2} & x_{4}^{U-2}-x_{U}x_{4}^{U-3} & x_{4}%
^{U-3}-x_{U}x_{4}^{U-4} & \cdots & x_{4}-x_{U} & 1\\
\vdots & \vdots & \vdots & \ddots & \vdots & \vdots\\
x_{U}^{U}-x_{U}^{2}x_{U}^{U-2} & x_{U}^{U-2}-x_{U}x_{U}^{U-3} & x_{U}%
^{U-3}-x_{U}x_{U}^{U-4} & \cdots & x_{U}-x_{U} & 1
\end{array}
\right)  .
\]
We claim that $\det B=\det\left(  A_{U}\right)  $. Indeed, here are two ways
to prove this:

\textit{First proof of }$\det B=\det\left(  A_{U}\right)  $\textit{:} Exercise
\ref{exe.ps4.6k} \textbf{(b)} shows that the determinant of a $U\times
U$-matrix does not change if we subtract a multiple of one of its columns from
another column. Now, let us do the following steps (in this order):

\begin{itemize}
\item subtract $x_{U}^{2}$ times the $2$-nd column of $A_{U}$ from the $1$-st column;

\item subtract $x_{U}$ times the $3$-rd column of the resulting matrix from
the $2$-nd column;

\item subtract $x_{U}$ times the $4$-th column of the resulting matrix from
the $3$-rd column;

\item and so on, all the way until we finally subtract $x_{U}$ times the
$U$-th column of the matrix from the $\left(  U-1\right)  $-st column.
\end{itemize}

Yes, you are reading this right: At the first step we subtract $x_{U}^{2}$
times (not $x_{U}$ times) the $2$-nd column from the $1$-st column; but at all
further steps, we subtract $x_{U}$ times a column from another. Having done
all this, the resulting matrix is $B$ (according to our definition of $B$).
Thus, $\det B=\det\left(  A_{U}\right)  $ (since our subtractions never change
the determinant). This proves $\det B=\det\left(  A_{U}\right)  $.

\textit{Second proof of }$\det B=\det\left(  A_{U}\right)  $\textit{:} Here is
another way to prove that $\det B=\det\left(  A_{U}\right)  $, with some less handwaving.

For every $\left(  i,j\right)  \in\left\{  1,2,\ldots,U\right\}  ^{2}$, we
define $c_{i,j}\in\mathbb{K}$ by%
\[
c_{i,j}=%
\begin{cases}
1, & \text{if }i=j;\\
-x_{U}^{2}, & \text{if }i=j+1\text{ and }j=1;\\
-x_{U}, & \text{if }i=j+1\text{ and }j>1;\\
0, & \text{otherwise}%
\end{cases}
.
\]
Let $C$ be the $U\times U$-matrix $\left(  c_{i,j}\right)  _{1\leq i\leq
U,\ 1\leq j\leq U}$. Here is how $C$ looks like:%
\[
C=\left(
\begin{array}
[c]{cccccc}%
1 & 0 & 0 & \cdots & 0 & 0\\
-x_{U}^{2} & 1 & 0 & \cdots & 0 & 0\\
0 & -x_{U} & 1 & \cdots & 0 & 0\\
0 & 0 & -x_{U} & \cdots & 0 & 0\\
\vdots & \vdots & \vdots & \ddots & \vdots & \vdots\\
0 & 0 & 0 & \cdots & -x_{U} & 1
\end{array}
\right)  ,
\]
where the only $-x_{U}^{2}$ is in the $\left(  2,1\right)  $-th cell.

The matrix $C$ is lower-triangular, and thus Exercise \ref{exe.ps4.3} shows
that its determinant is $\det C=\underbrace{c_{1,1}}_{=1}\underbrace{c_{2,2}%
}_{=1}\cdots\underbrace{c_{U,U}}_{=1}=1$.

On the other hand, it is easy to see that $B=A_{U}C$ (check this!). Thus,
Theorem \ref{thm.det(AB)} yields $\det B=\det\left(  A_{U}\right)
\cdot\underbrace{\det C}_{=1}=\det\left(  A_{U}\right)  $. So we have proven
$\det B=\det\left(  A_{U}\right)  $ again.

Next, we observe that for every $j\in\left\{  1,2,\ldots,U-1\right\}  $, we
have%
\begin{align*}
b_{U,j}  &  =%
\begin{cases}
x_{U}^{U}-x_{U}^{2}x_{U}^{U-2}, & \text{if }j=1;\\
x_{U}^{U-j}-x_{U}x_{U}^{U-j-1}, & \text{if }1<j<U;\\
1, & \text{if }j=U
\end{cases}
\ \ \ \ \ \ \ \ \ \ \left(  \text{by the definition of }b_{U,j}\right) \\
&  =%
\begin{cases}
x_{U}^{U}-x_{U}^{2}x_{U}^{U-2}, & \text{if }j=1;\\
x_{U}^{U-j}-x_{U}x_{U}^{U-j-1}, & \text{if }1<j<U
\end{cases}
\\
&  \ \ \ \ \ \ \ \ \ \ \left(  \text{since }j<U\text{ (since }j\in\left\{
1,2,\ldots,U-1\right\}  \text{)}\right) \\
&  =%
\begin{cases}
x_{U}^{U}-x_{U}^{U}, & \text{if }j=1;\\
x_{U}^{U-j}-x_{U}^{U-j}, & \text{if }1<j<U
\end{cases}
=%
\begin{cases}
0, & \text{if }j=1;\\
0, & \text{if }1<j<U
\end{cases}
=0.
\end{align*}
Hence, Theorem \ref{thm.laplace.pre} (applied to $U$, $B$ and $b_{i,j}$
instead of $n$, $A$ and $a_{i,j}$) yields%
\begin{equation}
\det B=b_{U,U}\cdot\det\left(  \left(  b_{i,j}\right)  _{1\leq i\leq
U-1,\ 1\leq j\leq U-1}\right)  . \label{sol.vander-det.s1.detB=prod}%
\end{equation}
Let $B^{\prime}$ denote the $\left(  U-1\right)  \times\left(  U-1\right)
$-matrix $\left(  b_{i,j}\right)  _{1\leq i\leq U-1,\ 1\leq j\leq U-1}$.

The definition of $b_{U,U}$ yields%
\begin{align*}
b_{U,U}  &  =%
\begin{cases}
x_{U}^{U}-x_{U}^{2}x_{U}^{U-2}, & \text{if }U=1;\\
x_{U}^{U-U}-x_{U}x_{U}^{U-U-1}, & \text{if }1<U<U;\\
1, & \text{if }U=U
\end{cases}
\ \ \ \ \ \ \ \ \ \ \left(  \text{by the definition of }b_{U,U}\right) \\
&  =1\ \ \ \ \ \ \ \ \ \ \left(  \text{since }U=U\right)  .
\end{align*}
Thus, (\ref{sol.vander-det.s1.detB=prod}) becomes%
\[
\det B=\underbrace{b_{U,U}}_{=1}\cdot\det\left(  \underbrace{\left(
b_{i,j}\right)  _{1\leq i\leq U-1,\ 1\leq j\leq U-1}}_{=B^{\prime}}\right)
=\det\left(  B^{\prime}\right)  .
\]
Compared with $\det B=\det\left(  A_{U}\right)  $, this yields%
\begin{equation}
\det\left(  A_{U}\right)  =\det\left(  B^{\prime}\right)  .
\label{sol.vander-det.s1.detAU=detB'}%
\end{equation}

Now, let us take a closer look at $B^{\prime}$. Indeed, every $\left(
i,j\right)  \in\left\{  1,2,\ldots,U-1\right\}  ^{2}$ satisfies%
\begin{align}
b_{i,j}  &  =%
\begin{cases}
x_{i}^{U}-x_{U}^{2}x_{i}^{U-2}, & \text{if }j=1;\\
x_{i}^{U-j}-x_{U}x_{i}^{U-j-1}, & \text{if }1<j<U;\\
1, & \text{if }j=U
\end{cases}
\ \ \ \ \ \ \ \ \ \ \left(  \text{by the definition of }b_{i,j}\right)
\nonumber\\
&  =%
\begin{cases}
x_{i}^{U}-x_{U}^{2}x_{i}^{U-2}, & \text{if }j=1;\\
x_{i}^{U-j}-x_{U}x_{i}^{U-j-1}, & \text{if }1<j<U
\end{cases}
\nonumber\\
&  \ \ \ \ \ \ \ \ \ \ \left(
\begin{array}
[c]{c}%
\text{since }j<U\text{ (since }j\in\left\{  1,2,\ldots,U-1\right\} \\
\text{(since }\left(  i,j\right)  \in\left\{  1,2,\ldots,U-1\right\}
^{2}\text{))}%
\end{array}
\right) \nonumber\\
&  =%
\begin{cases}
\left(  x_{i}^{2}-x_{U}^{2}\right)  x_{i}^{U-2}, & \text{if }j=1;\\
\left(  x_{i}-x_{U}\right)  x_{i}^{\left(  U-1\right)  -j}, & \text{if }1<j<U
\end{cases}
=%
\begin{cases}
\left(  x_{i}-x_{U}\right)  \left(  x_{i}+x_{U}\right)  x_{i}^{U-2}, &
\text{if }j=1;\\
\left(  x_{i}-x_{U}\right)  x_{i}^{\left(  U-1\right)  -j}, & \text{if }1<j<U
\end{cases}
\nonumber\\
&  =\left(  x_{i}-x_{U}\right)
\begin{cases}
\left(  x_{i}+x_{U}\right)  x_{i}^{U-2}, & \text{if }j=1;\\
x_{i}^{\left(  U-1\right)  -j}, & \text{if }1<j<U
\end{cases}
\nonumber\\
&  =\left(  x_{i}-x_{U}\right)
\begin{cases}
\left(  x_{i}+x_{U}\right)  x_{i}^{U-2}, & \text{if }j=1;\\
x_{i}^{\left(  U-1\right)  -j}, & \text{if }j>1
\end{cases}
\label{sol.vander-det.s1.bij-small}%
\end{align}
(since $1<j<U$ is equivalent to $j>1$).

For every $\left(  i,j\right)  \in\left\{  1,2,\ldots,U-1\right\}  ^{2}$, we
define $g_{i,j}\in\mathbb{K}$ by%
\[
g_{i,j}=%
\begin{cases}
\left(  x_{i}+x_{U}\right)  x_{i}^{U-2}, & \text{if }j=1;\\
x_{i}^{\left(  U-1\right)  -j}, & \text{if }j>1
\end{cases}
.
\]
Let $G$ be the $\left(  U-1\right)  \times\left(  U-1\right)  $-matrix
$\left(  g_{i,j}\right)  _{1\leq i\leq U-1,\ 1\leq j\leq U-1}$. Here is how
$G$ looks like:%
\[
G=\left(
\begin{array}
[c]{cccccc}%
\left(  x_{1}+x_{U}\right)  x_{1}^{U-2} & x_{1}^{U-3} & x_{1}^{U-4} & \cdots &
x_{1} & 1\\
\left(  x_{2}+x_{U}\right)  x_{2}^{U-2} & x_{2}^{U-3} & x_{2}^{U-4} & \cdots &
x_{2} & 1\\
\left(  x_{3}+x_{U}\right)  x_{3}^{U-2} & x_{3}^{U-3} & x_{3}^{U-4} & \cdots &
x_{3} & 1\\
\left(  x_{4}+x_{U}\right)  x_{4}^{U-2} & x_{4}^{U-3} & x_{4}^{U-4} & \cdots &
x_{4} & 1\\
\vdots & \vdots & \vdots & \ddots & \vdots & \vdots\\
\left(  x_{U-1}+x_{U}\right)  x_{U-1}^{U-2} & x_{U-1}^{U-3} & x_{U-1}^{U-4} &
\cdots & x_{U-1} & 1
\end{array}
\right)  .
\]
Now, (\ref{sol.vander-det.s1.bij-small}) becomes%
\begin{equation}
b_{i,j}=\left(  x_{i}-x_{U}\right)  \underbrace{%
\begin{cases}
\left(  x_{i}+x_{U}\right)  x_{i}^{U-2}, & \text{if }j=1;\\
x_{i}^{\left(  U-1\right)  -j}, & \text{if }j>1
\end{cases}
}_{=g_{i,j}}=\left(  x_{i}-x_{U}\right)  g_{i,j}
\label{sol.vander-det.s1.bij-small-2}%
\end{equation}
for every $\left(  i,j\right)  \in\left\{  1,2,\ldots,U-1\right\}  ^{2}$.
Hence,%
\begin{equation}
B^{\prime}=\left(  \underbrace{b_{i,j}}_{\substack{=\left(  x_{i}%
-x_{U}\right)  g_{i,j}\\\text{(by (\ref{sol.vander-det.s1.bij-small-2}))}%
}}\right)  _{1\leq i\leq U-1,\ 1\leq j\leq U-1}=\left(  \left(  x_{i}%
-x_{U}\right)  g_{i,j}\right)  _{1\leq i\leq U-1,\ 1\leq j\leq U-1}.
\label{sol.vander-det.s1.B'}%
\end{equation}
On the other hand, the definition of $G$ yields
\begin{equation}
G=\left(  g_{i,j}\right)  _{1\leq i\leq U-1,\ 1\leq j\leq U-1}.
\label{sol.vander-det.s1.AU-1}%
\end{equation}

Now, we claim that%
\begin{equation}
\det\left(  B^{\prime}\right)  =\det G\cdot\prod_{i=1}^{U-1}\left(
x_{i}-x_{U}\right)  . \label{sol.vander-det.s1.detB'=}%
\end{equation}
This can be proven similarly to how we proved (\ref{pf.thm.vander-det.detB'=})
back in our First proof of Theorem \ref{thm.vander-det}. (Of course, this
time, $G$ plays the role of $A_{U-1}$.) Now,
(\ref{sol.vander-det.s1.detAU=detB'}) becomes%
\begin{equation}
\det\left(  A_{U}\right)  =\det\left(  B^{\prime}\right)  =\det G\cdot
\prod_{i=1}^{U-1}\left(  x_{i}-x_{U}\right)  .
\label{sol.vander-det.s1.detAU=}%
\end{equation}

Now, we want to compute $\det G$. (This part is harder than the analogous part
of our First proof of Theorem \ref{thm.vander-det}, because back then the role
of $G$ was played by $A_{U-1}$, and we knew $\det\left(  A_{U-1}\right)  $
directly from our induction hypothesis.)

For every $\left(  i,j\right)  \in\left\{  1,2,\ldots,U-1\right\}  ^{2}$, we
define two elements $p_{i,j}\in\mathbb{K}$ and $q_{i,j}\in\mathbb{K}$ by%
\[
p_{i,j}=%
\begin{cases}
x_{i}^{U-1}, & \text{if }j=1;\\
x_{i}^{\left(  U-1\right)  -j}, & \text{if }j>1
\end{cases}
\]
and%
\[
q_{i,j}=%
\begin{cases}
x_{U}x_{i}^{U-2}, & \text{if }j=1;\\
x_{i}^{\left(  U-1\right)  -j}, & \text{if }j>1
\end{cases}
.
\]
We let $P$ be the $\left(  U-1\right)  \times\left(  U-1\right)  $-matrix
$\left(  p_{i,j}\right)  _{1\leq i\leq U-1,\ 1\leq j\leq U-1}$, and we let $Q$
be the $\left(  U-1\right)  \times\left(  U-1\right)  $-matrix $\left(
q_{i,j}\right)  _{1\leq i\leq U-1,\ 1\leq j\leq U-1}$. We now make the
following observations:

\begin{itemize}
\item The columns of the matrix $Q$ equal the corresponding columns of $P$
except (perhaps) the $1$-st column. The matrix $G$ is the $\left(  U-1\right)
\times\left(  U-1\right)  $-matrix obtained from $P$ by adding the $1$-st
column of $Q$ to the $1$-st column of $P$. Thus, Exercise \ref{exe.ps4.6}
\textbf{(j)} (applied to $U-1$, $1$, $P$, $Q$ and $G$ instead of $n$, $k$,
$A$, $A^{\prime}$ and $B$) yields $\det G=\det P+\det Q$.

\item We have $P=\left(  p_{i,j}\right)  _{1\leq i\leq U-1,\ 1\leq j\leq U-1}$
and $A_{U-1}=\left(  a_{i,j,U-1}\right)  _{1\leq i\leq U-1,\ 1\leq j\leq U-1}%
$. But a quick look at the definitions reveals that
\[
p_{i,j}=%
\begin{cases}
x_{i}^{U-1}, & \text{if }j=1;\\
x_{i}^{\left(  U-1\right)  -j}, & \text{if }j>1
\end{cases}
=%
\begin{cases}
x_{i}^{\left(  U-1\right)  -j}, & \text{if }j>1;\\
x_{i}^{U-1}, & \text{if }j=1
\end{cases}
=a_{i,j,U-1}%
\]
for all $\left(  i,j\right)  \in\left\{  1,2,\ldots,U-1\right\}  ^{2}$. Hence,
$P=A_{U-1}$.

\item The matrix $Q$ is obtained from the matrix $\left(  x_{i}^{\left(
U-1\right)  -j}\right)  _{1\leq i\leq U-1,\ 1\leq j\leq U-1}$ by multiplying
its $1$-st column by $x_{U}$. Hence, Exercise \ref{exe.ps4.6} \textbf{(h)}
(applied to $U-1$, $x_{U}$, $1$, $\left(  x_{i}^{\left(  U-1\right)
-j}\right)  _{1\leq i\leq U-1,\ 1\leq j\leq U-1}$ and $Q$ instead of $n$,
$\lambda$, $k$, $A$ and $B$) yields%
\begin{align*}
\det Q  &  =x_{U}\underbrace{\det\left(  \left(  x_{i}^{\left(  U-1\right)
-j}\right)  _{1\leq i\leq U-1,\ 1\leq j\leq U-1}\right)  }_{\substack{=\prod
_{1\leq i<j\leq U-1}\left(  x_{i}-x_{j}\right)  \\\text{(by Theorem
\ref{thm.vander-det} \textbf{(a)}, applied to }n=U-1\text{)}}}\\
&  =x_{U}\prod_{1\leq i<j\leq U-1}\left(  x_{i}-x_{j}\right)  .
\end{align*}

\end{itemize}

Combining this, we obtain%
\begin{align*}
\det G  &  =\det\underbrace{P}_{=A_{U-1}}+\underbrace{\det Q}_{=x_{U}%
\prod_{1\leq i<j\leq U-1}\left(  x_{i}-x_{j}\right)  }\\
&  =\underbrace{\det\left(  A_{U-1}\right)  }_{\substack{=\left(  x_{1}%
+x_{2}+\cdots+x_{U-1}\right)  \prod_{1\leq i<j\leq U-1}\left(  x_{i}%
-x_{j}\right)  \\\text{(since we assumed that (\ref{sol.vander-det.s1.goal})
holds for }u=U-1\text{)}}}+x_{U}\prod_{1\leq i<j\leq U-1}\left(  x_{i}%
-x_{j}\right) \\
&  =\left(  x_{1}+x_{2}+\cdots+x_{U-1}\right)  \prod_{1\leq i<j\leq
U-1}\left(  x_{i}-x_{j}\right)  +x_{U}\prod_{1\leq i<j\leq U-1}\left(
x_{i}-x_{j}\right) \\
&  =\underbrace{\left(  \left(  x_{1}+x_{2}+\cdots+x_{U-1}\right)
+x_{U}\right)  }_{=x_{1}+x_{2}+\cdots+x_{U}}\prod_{1\leq i<j\leq U-1}\left(
x_{i}-x_{j}\right) \\
&  =\left(  x_{1}+x_{2}+\cdots+x_{U}\right)  \prod_{1\leq i<j\leq U-1}\left(
x_{i}-x_{j}\right)  .
\end{align*}
Hence, (\ref{sol.vander-det.s1.detAU=}) yields%
\begin{align*}
\det\left(  A_{U}\right)   &  =\underbrace{\det G}_{=\left(  x_{1}%
+x_{2}+\cdots+x_{U}\right)  \prod_{1\leq i<j\leq U-1}\left(  x_{i}%
-x_{j}\right)  }\cdot\prod_{i=1}^{U-1}\left(  x_{i}-x_{U}\right) \\
&  =\left(  x_{1}+x_{2}+\cdots+x_{U}\right)  \underbrace{\prod_{1\leq i<j\leq
U-1}}_{=\prod_{j=1}^{U-1}\prod_{i=1}^{j-1}}\left(  x_{i}-x_{j}\right)
\cdot\prod_{i=1}^{U-1}\left(  x_{i}-x_{U}\right) \\
&  =\left(  x_{1}+x_{2}+\cdots+x_{U}\right)  \left(  \prod_{j=1}^{U-1}%
\prod_{i=1}^{j-1}\left(  x_{i}-x_{j}\right)  \right)  \cdot\prod_{i=1}%
^{U-1}\left(  x_{i}-x_{U}\right)  .
\end{align*}
Compared with%
\begin{align*}
&  \left(  x_{1}+x_{2}+\cdots+x_{U}\right)  \underbrace{\prod_{1\leq i<j\leq
U}}_{=\prod_{j=1}^{U}\prod_{i=1}^{j-1}}\left(  x_{i}-x_{j}\right) \\
&  =\left(  x_{1}+x_{2}+\cdots+x_{U}\right)  \prod_{j=1}^{U}\prod_{i=1}%
^{j-1}\left(  x_{i}-x_{j}\right) \\
&  =\left(  x_{1}+x_{2}+\cdots+x_{U}\right)  \left(  \prod_{j=1}^{U-1}%
\prod_{i=1}^{j-1}\left(  x_{i}-x_{j}\right)  \right)  \cdot\prod_{i=1}%
^{U-1}\left(  x_{i}-x_{U}\right) \\
&  \ \ \ \ \ \ \ \ \ \ \left(  \text{here, we have split off the factor for
}j=U\text{ from the product}\right)  ,
\end{align*}
this yields $\det\left(  A_{U}\right)  =\left(  x_{1}+x_{2}+\cdots
+x_{U}\right)  \prod_{1\leq i<j\leq U}\left(  x_{i}-x_{j}\right)  $. In other
words, (\ref{sol.vander-det.s1.goal}) holds for $u=U$. This completes the
induction step.

Now, (\ref{sol.vander-det.s1.goal}) is proven by induction.]

Hence, we can apply (\ref{sol.vander-det.s1.goal}) to $u=n$. As the result, we
obtain $\det\left(  A_{n}\right)  =\left(  x_{1}+x_{2}+\cdots+x_{n}\right)
\prod_{1\leq i<j\leq n}\left(  x_{i}-x_{j}\right)  $. Since $A_{n}=\left(
\begin{cases}
x_{i}^{n-j}, & \text{if }j>1;\\
x_{i}^{n}, & \text{if }j=1
\end{cases}
\right)  _{1\leq i\leq n,\ 1\leq j\leq n}$ (by the definition of $A_{n}$),
this rewrites as
\[
\det\left(  \left(
\begin{cases}
x_{i}^{n-j}, & \text{if }j>1;\\
x_{i}^{n}, & \text{if }j=1
\end{cases}
\right)  _{1\leq i\leq n,\ 1\leq j\leq n}\right)  =\left(  x_{1}+x_{2}%
+\cdots+x_{n}\right)  \prod_{1\leq i<j\leq n}\left(  x_{i}-x_{j}\right)  .
\]
This solves Exercise \ref{exe.vander-det.s1}.
\end{proof}

We shall give a second solution for Exercise \ref{exe.vander-det.s1} later (in
Section \ref{sect.sols.vander-det.s1.sol2}).

\subsection{Solution to Exercise \ref{exe.vander-det.xi+yj}}

We shall first sketch short solutions to parts \textbf{(a)} and \textbf{(b)}
of Exercise \ref{exe.vander-det.xi+yj}. Then, we will show a solution to part
\textbf{(c)} (which is a more elaborate and subtler version of our solution to
part \textbf{(b)}) and a solution to part \textbf{(d)} (which is quite similar
to that to part \textbf{(c)}), and finally derive parts \textbf{(a)} and
\textbf{(b)} from part \textbf{(c)}. Thus, parts \textbf{(a)} and \textbf{(b)}
of Exercise \ref{exe.vander-det.xi+yj} will be solved twice (though the
solutions cannot really be called different).

\begin{proof}
[Solution sketch to parts \textbf{(a)} and \textbf{(b)} of Exercise
\ref{exe.vander-det.xi+yj}.]\textbf{(a)} Let $m\in\left\{  0,1,\ldots
,n-2\right\}  $. Thus, $m+1<n$.

Define an $n\times\left(  m+1\right)  $-matrix $B$ by%
\[
B=\left(  \dbinom{m}{j-1}x_{i}^{j-1}\right)  _{1\leq i\leq n,\ 1\leq j\leq
m+1}.
\]
Define an $\left(  m+1\right)  \times n$-matrix $C$ by%
\[
C=\left(  y_{j}^{m-\left(  i-1\right)  }\right)  _{1\leq i\leq m+1,\ 1\leq
j\leq n}.
\]
Then, the definition of the product of two matrices shows that%
\[
BC=\left(  \sum_{k=1}^{m+1}\dbinom{m}{k-1}x_{i}^{k-1}y_{j}^{m-\left(
k-1\right)  }\right)  _{1\leq i\leq n,\ 1\leq j\leq n}.
\]
Since every $\left(  i,j\right)  \in\left\{  1,2,\ldots,n\right\}  ^{2}$
satisfies%
\begin{align*}
&  \sum_{k=1}^{m+1}\dbinom{m}{k-1}x_{i}^{k-1}y_{j}^{m-\left(  k-1\right)  }\\
&  =\sum_{k=0}^{m}\dbinom{m}{k}x_{i}^{k}y_{j}^{m-k}\ \ \ \ \ \ \ \ \ \ \left(
\text{here, we have substituted }k\text{ for }k-1\text{ in the sum}\right) \\
&  =\left(  x_{i}+y_{j}\right)  ^{m}\ \ \ \ \ \ \ \ \ \ \left(
\begin{array}
[c]{c}%
\text{since the binomial formula yields}\\
\left(  x_{i}+y_{j}\right)  ^{m}=\sum_{k=0}^{m}\dbinom{m}{k}x_{i}^{k}%
y_{j}^{m-k}%
\end{array}
\right)  ,
\end{align*}
this rewrites as%
\[
BC=\left(  \left(  x_{i}+y_{j}\right)  ^{m}\right)  _{1\leq i\leq n,\ 1\leq
j\leq n}.
\]

But recall that $m+1<n$. Hence, (\ref{eq.exam.cauchy-binet.0}) (applied to
$m+1$, $B$ and $C$ instead of $m$, $A$ and $B$) shows that $\det\left(
BC\right)  =0$. Since $BC=\left(  \left(  x_{i}+y_{j}\right)  ^{m}\right)
_{1\leq i\leq n,\ 1\leq j\leq n}$, this rewrites as $\det\left(  \left(
\left(  x_{i}+y_{j}\right)  ^{m}\right)  _{1\leq i\leq n,\ 1\leq j\leq
n}\right)  =0$. This solves Exercise \ref{exe.vander-det.xi+yj} \textbf{(a)}.

\textbf{(b)} Define an $n\times n$-matrix $B$ by%
\[
B=\left(  \dbinom{n-1}{j-1}x_{i}^{j-1}\right)  _{1\leq i\leq n,\ 1\leq j\leq
n}.
\]
Define an $n\times n$-matrix $C$ by%
\[
C=\left(  y_{j}^{n-i}\right)  _{1\leq i\leq n,\ 1\leq j\leq n}.
\]
Then, the definition of the product of two matrices shows that%
\[
BC=\left(  \sum_{k=1}^{n}\dbinom{n-1}{k-1}x_{i}^{k-1}y_{j}^{n-k}\right)
_{1\leq i\leq n,\ 1\leq j\leq n}.
\]
Since every $\left(  i,j\right)  \in\left\{  1,2,\ldots,n\right\}  ^{2}$
satisfies%
\begin{align*}
&  \sum_{k=1}^{n}\dbinom{n-1}{k-1}x_{i}^{k-1}\underbrace{y_{j}^{n-k}}%
_{=y_{j}^{\left(  n-1\right)  -\left(  k-1\right)  }}=\sum_{k=1}^{n}%
\dbinom{n-1}{k-1}x_{i}^{k-1}y_{j}^{\left(  n-1\right)  -\left(  k-1\right)
}\\
&  =\sum_{k=0}^{n-1}\dbinom{n-1}{k}x_{i}^{k}y_{j}^{\left(  n-1\right)
-k}\ \ \ \ \ \ \ \ \ \ \left(  \text{here, we have substituted }k\text{ for
}k-1\text{ in the sum}\right) \\
&  =\left(  x_{i}+y_{j}\right)  ^{n-1}\ \ \ \ \ \ \ \ \ \ \left(
\begin{array}
[c]{c}%
\text{since the binomial formula yields}\\
\left(  x_{i}+y_{j}\right)  ^{n-1}=\sum_{k=0}^{n-1}\dbinom{n-1}{k}x_{i}%
^{k}y_{j}^{\left(  n-1\right)  -k}%
\end{array}
\right)  ,
\end{align*}
this rewrites as%
\[
BC=\left(  \left(  x_{i}+y_{j}\right)  ^{n-1}\right)  _{1\leq i\leq n,\ 1\leq
j\leq n}.
\]

Theorem \ref{thm.det(AB)} shows that $\det\left(  BC\right)  =\det B\cdot\det
C$. But what are $\det B$ and $\det C$ ?

Finding $\det C$ is easy: We have $C=\left(  y_{j}^{n-i}\right)  _{1\leq i\leq
n,\ 1\leq j\leq n}$ and thus
\[
\det C=\det\left(  \left(  y_{j}^{n-i}\right)  _{1\leq i\leq n,\ 1\leq j\leq
n}\right)  =\prod_{1\leq i<j\leq n}\left(  y_{i}-y_{j}\right)
\]
(by Theorem \ref{thm.vander-det} \textbf{(b)}, applied to $y_{k}$ instead of
$x_{k}$).

To find $\det B$, we first observe that $\det\left(  \left(  x_{i}%
^{j-1}\right)  _{1\leq i\leq n,\ 1\leq j\leq n}\right)  =\prod_{1\leq j<i\leq
n}\left(  x_{i}-x_{j}\right)  $ (by Theorem \ref{thm.vander-det}
\textbf{(c)}). But $B=\left(  \dbinom{n-1}{j-1}x_{i}^{j-1}\right)  _{1\leq
i\leq n,\ 1\leq j\leq n}$. In other words, the matrix $B$ is obtained from the
matrix $\left(  x_{i}^{j-1}\right)  _{1\leq i\leq n,\ 1\leq j\leq n}$ by
multiplying the whole $j$-th column with $\dbinom{n-1}{j-1}$ for every
$j\in\left\{  1,2,\ldots,n\right\}  $. Therefore, the determinant $\det B$ is
obtained by successively multiplying $\det\left(  \left(  x_{i}^{j-1}\right)
_{1\leq i\leq n,\ 1\leq j\leq n}\right)  $ with $\dbinom{n-1}{j-1}$ for every
$j\in\left\{  1,2,\ldots,n\right\}  $ (since Exercise \ref{exe.ps4.6}
\textbf{(h)} tells us that multiplying a single column of a square matrix by a
scalar $\lambda$ results in the determinant getting multiplied by $\lambda$).
In other words,%
\begin{align*}
\det B  &  =\underbrace{\det\left(  \left(  x_{i}^{j-1}\right)  _{1\leq i\leq
n,\ 1\leq j\leq n}\right)  }_{\substack{=\prod_{1\leq j<i\leq n}\left(
x_{i}-x_{j}\right)  \\=\prod_{1\leq i<j\leq n}\left(  x_{j}-x_{i}\right)
\\\text{(here, we renamed the}\\\text{index }\left(  j,i\right)  \text{ as
}\left(  i,j\right)  \text{)}}}\cdot\underbrace{\prod_{j=1}^{n}\dbinom
{n-1}{j-1}}_{\substack{=\prod_{k=0}^{n-1}\dbinom{n-1}{k}\\\text{(here, we
substituted }k\text{ for }j-1\text{)}}}\\
&  =\left(  \prod_{1\leq i<j\leq n}\left(  x_{j}-x_{i}\right)  \right)
\cdot\left(  \prod_{k=0}^{n-1}\dbinom{n-1}{k}\right)  .
\end{align*}

Now,%
\begin{align*}
\det\left(  BC\right)   &  =\underbrace{\det B}_{=\left(  \prod_{1\leq i<j\leq
n}\left(  x_{j}-x_{i}\right)  \right)  \cdot\left(  \prod_{k=0}^{n-1}%
\dbinom{n-1}{k}\right)  }\cdot\underbrace{\det C}_{=\prod_{1\leq i<j\leq
n}\left(  y_{i}-y_{j}\right)  }\\
&  =\left(  \prod_{1\leq i<j\leq n}\left(  x_{j}-x_{i}\right)  \right)
\cdot\left(  \prod_{k=0}^{n-1}\dbinom{n-1}{k}\right)  \cdot\left(
\prod_{1\leq i<j\leq n}\left(  y_{i}-y_{j}\right)  \right) \\
&  =\left(  \prod_{k=0}^{n-1}\dbinom{n-1}{k}\right)  \cdot\underbrace{\left(
\prod_{1\leq i<j\leq n}\left(  x_{j}-x_{i}\right)  \right)  \cdot\left(
\prod_{1\leq i<j\leq n}\left(  y_{i}-y_{j}\right)  \right)  }_{=\prod_{1\leq
i<j\leq n}\left(  \left(  x_{j}-x_{i}\right)  \left(  y_{i}-y_{j}\right)
\right)  }\\
&  =\left(  \prod_{k=0}^{n-1}\dbinom{n-1}{k}\right)  \cdot\prod_{1\leq i<j\leq
n}\underbrace{\left(  \left(  x_{j}-x_{i}\right)  \left(  y_{i}-y_{j}\right)
\right)  }_{=\left(  x_{i}-x_{j}\right)  \left(  y_{j}-y_{i}\right)  }\\
&  =\left(  \prod_{k=0}^{n-1}\dbinom{n-1}{k}\right)  \cdot\underbrace{\prod
_{1\leq i<j\leq n}\left(  \left(  x_{i}-x_{j}\right)  \left(  y_{j}%
-y_{i}\right)  \right)  }_{=\left(  \prod_{1\leq i<j\leq n}\left(  x_{i}%
-x_{j}\right)  \right)  \left(  \prod_{1\leq i<j\leq n}\left(  y_{j}%
-y_{i}\right)  \right)  }\\
&  =\left(  \prod_{k=0}^{n-1}\dbinom{n-1}{k}\right)  \cdot\left(  \prod_{1\leq
i<j\leq n}\left(  x_{i}-x_{j}\right)  \right)  \left(  \prod_{1\leq i<j\leq
n}\left(  y_{j}-y_{i}\right)  \right)  .
\end{align*}
Since $BC=\left(  \left(  x_{i}+y_{j}\right)  ^{n-1}\right)  _{1\leq i\leq
n,\ 1\leq j\leq n}$, this rewrites as%
\begin{align*}
&  \det\left(  \left(  \left(  x_{i}+y_{j}\right)  ^{n-1}\right)  _{1\leq
i\leq n,\ 1\leq j\leq n}\right) \\
&  =\left(  \prod_{k=0}^{n-1}\dbinom{n-1}{k}\right)  \cdot\left(  \prod_{1\leq
i<j\leq n}\left(  x_{i}-x_{j}\right)  \right)  \left(  \prod_{1\leq i<j\leq
n}\left(  y_{j}-y_{i}\right)  \right)  .
\end{align*}
This solves Exercise \ref{exe.vander-det.xi+yj} \textbf{(b)}.
\end{proof}

Now, as promised, we shall start from scratch, and solve Exercise
\ref{exe.vander-det.xi+yj} \textbf{(c)} and \textbf{(d)} first, and then solve
Exercise \ref{exe.vander-det.xi+yj} \textbf{(a)} and \textbf{(b)} again using
Exercise \ref{exe.vander-det.xi+yj} \textbf{(c)}.

\begin{vershort}
\begin{proof}
[Solution to Exercise \ref{exe.vander-det.xi+yj}.]\textbf{(c)} We extend the
$n$-tuple $\left(  p_{0},p_{1},\ldots,p_{n-1}\right)  $ to an infinite
sequence $\left(  p_{0},p_{1},p_{2},\ldots\right)  \in\mathbb{K}^{\infty}$ by
setting%
\begin{equation}
p_{\ell}=0\ \ \ \ \ \ \ \ \ \ \text{for every }\ell\in\left\{
n,n+1,n+2,\ldots\right\}  . \label{sol.vander-det.xi+yj.short.c.pk=0}%
\end{equation}

Next, we define three $n\times n$-matrices $B$, $C$ and $D$:

\begin{itemize}
\item Define an $n\times n$-matrix $B$ by
\[
B=\left(  x_{i}^{n-j}\right)  _{1\leq i\leq n,\ 1\leq j\leq n}.
\]
Then,%
\begin{align}
\det\underbrace{B}_{=\left(  x_{i}^{n-j}\right)  _{1\leq i\leq n,\ 1\leq j\leq
n}}  &  =\det\left(  \left(  x_{i}^{n-j}\right)  _{1\leq i\leq n,\ 1\leq j\leq
n}\right) \nonumber\\
&  =\prod_{1\leq i<j\leq n}\left(  x_{i}-x_{j}\right)
\label{sol.vander-det.xi+yj.short.c.detB}%
\end{align}
(by Theorem \ref{thm.vander-det} \textbf{(a)}).

\item For every $\left(  i,j\right)  \in\left\{  1,2,\ldots,n\right\}  ^{2}$,
the element $\dbinom{n-1-i+j}{j-1}p_{n-1-i+j}$ of $\mathbb{K}$ is
well-defined\footnote{\textit{Proof.} Let $\left(  i,j\right)  \in\left\{
1,2,\ldots,n\right\}  ^{2}$. Thus, $j-1\geq0$. Hence, the binomial coefficient
$\dbinom{n-1-i+j}{j-1}\in\mathbb{Z}$ is well-defined. Moreover,
$n-1-\underbrace{i}_{\leq n}+\underbrace{j}_{\geq1}\geq n-1-n+1=0$, so that
$n-1-i+j\in\mathbb{N}$. Thus, $p_{n-1-i+j}\in\mathbb{K}$ is well-defined
(since we have an infinite sequence $\left(  p_{0},p_{1},p_{2},\ldots\right)
\in\mathbb{K}^{\infty}$). Hence, the element $\dbinom{n-1-i+j}{j-1}%
p_{n-1-i+j}$ of $\mathbb{K}$ is well-defined. Qed.}. Thus, for every $\left(
i,j\right)  \in\left\{  1,2,\ldots,n\right\}  ^{2}$, we can define an element
$c_{i,j}\in\mathbb{K}$ by $c_{i,j}=\dbinom{n-1-i+j}{j-1}p_{n-1-i+j}$. Consider
these elements $c_{i,j}$. Define an $n\times n$-matrix $C$ by%
\[
C=\left(  c_{i,j}\right)  _{1\leq i\leq n,\ 1\leq j\leq n}.
\]
Consider this matrix $C$. We have $c_{i,j}=0$ for every $\left(  i,j\right)
\in\left\{  1,2,\ldots,n\right\}  ^{2}$ satisfying $i<j$%
\ \ \ \ \footnote{\textit{Proof.} Let $\left(  i,j\right)  \in\left\{
1,2,\ldots,n\right\}  ^{2}$ be such that $i<j$. From $i<j$, we obtain $i\leq
j-1$ (since $i$ and $j$ are integers). Thus, $n-1-\underbrace{i}_{\leq
j-1}+j\geq n-1-\left(  j-1\right)  +j=n$. In other words, $n-1-i+j\in\left\{
n,n+1,n+2,\ldots\right\}  $. Hence, $p_{n-1-i+j}=0$ (by
(\ref{sol.vander-det.xi+yj.short.c.pk=0}), applied to $\ell=n-1-i+j$). Hence,
$c_{i,j}=\dbinom{n-1-i+j}{j-1}\underbrace{p_{n-1-i+j}}_{=0}=0$, qed.}. Hence,
Exercise \ref{exe.ps4.3} (applied to $C$ and $\left(  c_{i,j}\right)  _{1\leq
i\leq n,\ 1\leq j\leq n}$ instead of $A$ and $\left(  a_{i,j}\right)  _{1\leq
i\leq n,\ 1\leq j\leq n}$) shows that%
\begin{align}
\det C  &  =c_{1,1}c_{2,2}\cdots c_{n,n}=\prod_{i=1}^{n}\underbrace{c_{i,i}%
}_{\substack{=\dbinom{n-1-i+i}{i-1}p_{n-1-i+i}\\\text{(by the definition of
}c_{i,i}\text{)}}}\nonumber\\
&  =\prod_{i=1}^{n}\left(  \underbrace{\dbinom{n-1-i+i}{i-1}}_{=\dbinom
{n-1}{i-1}}\underbrace{p_{n-1-i+i}}_{=p_{n-1}}\right) \nonumber\\
&  =\prod_{i=1}^{n}\left(  \dbinom{n-1}{i-1}p_{n-1}\right)  =\left(
\prod_{i=1}^{n}\dbinom{n-1}{i-1}\right)  p_{n-1}^{n}\nonumber\\
&  =p_{n-1}^{n}\left(  \prod_{i=1}^{n}\dbinom{n-1}{i-1}\right)  =p_{n-1}%
^{n}\left(  \prod_{k=0}^{n-1}\dbinom{n-1}{k}\right)
\label{sol.vander-det.xi+yj.short.c.detC}%
\end{align}
(here, we have substituted $k$ for $i-1$ in the product).

\item Define an $n\times n$-matrix $D$ by
\[
D=\left(  y_{j}^{i-1}\right)  _{1\leq i\leq n,\ 1\leq j\leq n}.
\]
Then,%
\begin{align}
\det\underbrace{D}_{=\left(  y_{j}^{i-1}\right)  _{1\leq i\leq n,\ 1\leq j\leq
n}}  &  =\det\left(  \left(  y_{j}^{i-1}\right)  _{1\leq i\leq n,\ 1\leq j\leq
n}\right)  =\prod_{1\leq j<i\leq n}\left(  y_{i}-y_{j}\right) \nonumber\\
&  \ \ \ \ \ \ \ \ \ \ \left(  \text{by Theorem \ref{thm.vander-det}
\textbf{(d)}, applied to }y_{k}\text{ instead of }x_{k}\right) \nonumber\\
&  =\prod_{1\leq i<j\leq n}\left(  y_{j}-y_{i}\right)
\label{sol.vander-det.xi+yj.short.c.detD}%
\end{align}
(here, we have renamed the index $\left(  j,i\right)  $ as $\left(
i,j\right)  $ in the product).
\end{itemize}

Our goal is to prove that $\left(  P\left(  x_{i}+y_{j}\right)  \right)
_{1\leq i\leq n,\ 1\leq j\leq n}=BCD$. Once this is done, we will be able to
compute $\det\left(  \left(  P\left(  x_{i}+y_{j}\right)  \right)  _{1\leq
i\leq n,\ 1\leq j\leq n}\right)  $ by applying Theorem \ref{thm.det(AB)} twice.

We have $C=\left(  c_{i,j}\right)  _{1\leq i\leq n,\ 1\leq j\leq n}$ and
$D=\left(  y_{j}^{i-1}\right)  _{1\leq i\leq n,\ 1\leq j\leq n}$. Hence, the
definition of the product of two matrices shows that%
\begin{equation}
CD=\left(  \sum_{k=1}^{n}c_{i,k}y_{j}^{k-1}\right)  _{1\leq i\leq n,\ 1\leq
j\leq n}. \label{sol.vander-det.xi+yj.short.c.CD.1}%
\end{equation}
But every $\left(  i,j\right)  \in\left\{  1,2,\ldots,n\right\}  ^{2}$
satisfies%
\begin{align*}
\sum_{k=1}^{n}c_{i,k}y_{j}^{k-1}  &  =\sum_{\ell=0}^{n-1}\underbrace{c_{i,\ell
+1}}_{\substack{=\dbinom{n-1-i+\left(  \ell+1\right)  }{\left(  \ell+1\right)
-1}p_{n-1-i+\left(  \ell+1\right)  }\\\text{(by the definition of }%
c_{i,\ell+1}\text{)}}}y_{j}^{\left(  \ell+1\right)  -1}\\
&  \ \ \ \ \ \ \ \ \ \ \left(  \text{here, we have substituted }\ell+1\text{
for }k\text{ in the sum}\right) \\
&  =\sum_{\ell=0}^{n-1}\underbrace{\dbinom{n-1-i+\left(  \ell+1\right)
}{\left(  \ell+1\right)  -1}}_{\substack{=\dbinom{n-i+\ell}{\ell}%
}}\underbrace{p_{n-1-i+\left(  \ell+1\right)  }}_{=p_{n-i+\ell}}%
\underbrace{y_{j}^{\left(  \ell+1\right)  -1}}_{=y_{j}^{\ell}}\\
&  =\sum_{\ell=0}^{n-1}\dbinom{n-i+\ell}{\ell}p_{n-i+\ell}y_{j}^{\ell}.
\end{align*}
Hence, (\ref{sol.vander-det.xi+yj.short.c.CD.1}) rewrites as%
\[
CD=\left(  \sum_{\ell=0}^{n-1}\dbinom{n-i+\ell}{\ell}p_{n-i+\ell}y_{j}^{\ell
}\right)  _{1\leq i\leq n,\ 1\leq j\leq n}.
\]

Now, $B=\left(  x_{i}^{n-j}\right)  _{1\leq i\leq n,\ 1\leq j\leq n}$ and
$CD=\left(  \sum_{\ell=0}^{n-1}\dbinom{n-i+\ell}{\ell}p_{n-i+\ell}y_{j}^{\ell
}\right)  _{1\leq i\leq n,\ 1\leq j\leq n}$. Hence, the definition of the
product of two matrices shows that%
\begin{equation}
B\left(  CD\right)  =\left(  \sum_{k=1}^{n}x_{i}^{n-k}\left(  \sum_{\ell
=0}^{n-1}\dbinom{n-k+\ell}{\ell}p_{n-k+\ell}y_{j}^{\ell}\right)  \right)
_{1\leq i\leq n,\ 1\leq j\leq n}. \label{sol.vander-det.xi+yj.short.c.BCD.1}%
\end{equation}

Now, let $\left(  i,j\right)  \in\left\{  1,2,\ldots,n\right\}  ^{2}$. Then,%
\begin{align}
&  \sum_{k=1}^{n}x_{i}^{n-k}\left(  \sum_{\ell=0}^{n-1}\dbinom{n-k+\ell}{\ell
}p_{n-k+\ell}y_{j}^{\ell}\right) \nonumber\\
&  =\sum_{k=0}^{n-1}x_{i}^{k}\left(  \sum_{\ell=0}^{n-1}\dbinom{k+\ell}{\ell
}p_{k+\ell}y_{j}^{\ell}\right) \nonumber\\
&  \ \ \ \ \ \ \ \ \ \ \left(  \text{here, we have substituted }k\text{ for
}n-k\text{ in the first sum}\right) \nonumber\\
&  =\sum_{k=0}^{n-1}\sum_{\ell=0}^{n-1}x_{i}^{k}\dbinom{k+\ell}{\ell}%
p_{k+\ell}y_{j}^{\ell}\nonumber\\
&  =\sum_{k=0}^{n-1}\sum_{\ell=k}^{n-1+k}x_{i}^{k}\underbrace{\dbinom
{k+\left(  \ell-k\right)  }{\ell-k}}_{=\dbinom{\ell}{\ell-k}}%
\underbrace{p_{k+\left(  \ell-k\right)  }}_{=p_{\ell}}y_{j}^{\ell
-k}\nonumber\\
&  \ \ \ \ \ \ \ \ \ \ \left(  \text{here, we have substituted }\ell-k\text{
for }\ell\text{ in the second sum}\right) \nonumber\\
&  =\sum_{k=0}^{n-1}\underbrace{\sum_{\ell=k}^{n-1+k}x_{i}^{k}\dbinom{\ell
}{\ell-k}p_{\ell}y_{j}^{\ell-k}}_{\substack{=\sum_{\ell=k}^{n-1}x_{i}%
^{k}\dbinom{\ell}{\ell-k}p_{\ell}y_{j}^{\ell-k}+\sum_{\ell=n}^{n-1+k}x_{i}%
^{k}\dbinom{\ell}{\ell-k}p_{\ell}y_{j}^{\ell-k}\\\text{(since }k\leq n\leq
n-1+k+1\text{ (since }n\leq n+k=n-1+k+1\text{))}}}\nonumber\\
&  =\sum_{k=0}^{n-1}\left(  \sum_{\ell=k}^{n-1}x_{i}^{k}\dbinom{\ell}{\ell
-k}p_{\ell}y_{j}^{\ell-k}+\sum_{\ell=n}^{n-1+k}x_{i}^{k}\dbinom{\ell}{\ell
-k}\underbrace{p_{\ell}}_{\substack{=0\\\text{(by
(\ref{sol.vander-det.xi+yj.short.c.pk=0})}\\\text{(since }\ell\in\left\{
n,n+1,\ldots,n-1+k\right\}  \\\subseteq\left\{  n,n+1,n+2,\ldots\right\}
\text{))}}}y_{j}^{\ell-k}\right) \nonumber\\
&  =\sum_{k=0}^{n-1}\left(  \sum_{\ell=k}^{n-1}x_{i}^{k}\dbinom{\ell}{\ell
-k}p_{\ell}y_{j}^{\ell-k}+\underbrace{\sum_{\ell=n}^{n-1+k}x_{i}^{k}%
\dbinom{\ell}{\ell-k}0y_{j}^{\ell-k}}_{=0}\right) \nonumber\\
&  =\sum_{k=0}^{n-1}\sum_{\ell=k}^{n-1}x_{i}^{k}\dbinom{\ell}{\ell-k}p_{\ell
}y_{j}^{\ell-k}. \label{sol.vander-det.xi+yj.short.c.BCD.2}%
\end{align}
On the other hand, $P\left(  X\right)  =\sum_{k=0}^{n-1}p_{k}X^{k}=\sum
_{\ell=0}^{n-1}p_{\ell}X^{\ell}$ (here, we renamed the summation index $k$ as
$\ell$). Hence,%
\begin{align}
P\left(  x_{i}+y_{j}\right)   &  =\sum_{\ell=0}^{n-1}p_{\ell}%
\underbrace{\left(  x_{i}+y_{j}\right)  ^{\ell}}_{\substack{=\sum_{k=0}^{\ell
}\dbinom{\ell}{k}x_{i}^{k}y_{j}^{\ell-k}\\\text{(by the binomial formula)}%
}}=\sum_{\ell=0}^{n-1}p_{\ell}\sum_{k=0}^{\ell}\dbinom{\ell}{k}x_{i}^{k}%
y_{j}^{\ell-k}\nonumber\\
&  =\underbrace{\sum_{\ell=0}^{n-1}\sum_{k=0}^{\ell}}_{=\sum
_{\substack{\left(  \ell,k\right)  \in\left\{  0,1,\ldots,n-1\right\}
^{2};\\k\leq\ell}}=\sum_{k=0}^{n-1}\sum_{\ell=k}^{n-1}}p_{\ell}\dbinom{\ell
}{k}x_{i}^{k}y_{j}^{\ell-k}\nonumber\\
&  =\sum_{k=0}^{n-1}\sum_{\ell=k}^{n-1}p_{\ell}\underbrace{\dbinom{\ell}{k}%
}_{\substack{=\dbinom{\ell}{\ell-k}\\\text{(by (\ref{eq.binom.symm}), applied
to }\ell\text{ and }k\\\text{instead of }m\text{ and }n\text{)}}}x_{i}%
^{k}y_{j}^{\ell-k}\nonumber\\
&  =\sum_{k=0}^{n-1}\sum_{\ell=k}^{n-1}p_{\ell}\dbinom{\ell}{\ell-k}x_{i}%
^{k}y_{j}^{\ell-k}=\sum_{k=0}^{n-1}\sum_{\ell=k}^{n-1}x_{i}^{k}\dbinom{\ell
}{\ell-k}p_{\ell}y_{j}^{\ell-k}\nonumber\\
&  =\sum_{k=1}^{n}x_{i}^{n-k}\left(  \sum_{\ell=0}^{n-1}\dbinom{n-k+\ell}%
{\ell}p_{n-k+\ell}y_{j}^{\ell}\right)  \ \ \ \ \ \ \ \ \ \ \left(  \text{by
(\ref{sol.vander-det.xi+yj.short.c.BCD.2})}\right)  .
\label{sol.vander-det.xi+yj.short.c.BCD.3}%
\end{align}

Now, let us forget that we fixed $\left(  i,j\right)  $. We thus have proven
(\ref{sol.vander-det.xi+yj.short.c.BCD.3}) for every $\left(  i,j\right)
\in\left\{  1,2,\ldots,n\right\}  ^{2}$. Hence,%
\begin{align*}
&  \left(  P\left(  x_{i}+y_{j}\right)  \right)  _{1\leq i\leq n,\ 1\leq j\leq
n}\\
&  =\left(  \sum_{k=1}^{n}x_{i}^{n-k}\left(  \sum_{\ell=0}^{n-1}%
\dbinom{n-k+\ell}{\ell}p_{n-k+\ell}y_{j}^{\ell}\right)  \right)  _{1\leq i\leq
n,\ 1\leq j\leq n}=B\left(  CD\right)  \ \ \ \ \ \ \ \ \ \ \left(  \text{by
(\ref{sol.vander-det.xi+yj.short.c.BCD.1})}\right) \\
&  =BCD.
\end{align*}
Hence,%
\begin{align*}
&  \det\underbrace{\left(  \left(  P\left(  x_{i}+y_{j}\right)  \right)
_{1\leq i\leq n,\ 1\leq j\leq n}\right)  }_{=BCD}\\
&  =\det\left(  BCD\right)  =\det B\cdot\underbrace{\det\left(  CD\right)
}_{\substack{=\det C\cdot\det D\\\text{(by Theorem \ref{thm.det(AB)}, applied
to }C\text{ and }D\\\text{instead of }A\text{ and }B\text{)}}}\\
&  \ \ \ \ \ \ \ \ \ \ \left(  \text{by Theorem \ref{thm.det(AB)}, applied to
}B\text{ and }CD\text{ instead of }A\text{ and }B\right) \\
&  =\underbrace{\det B}_{\substack{=\prod_{1\leq i<j\leq n}\left(  x_{i}%
-x_{j}\right)  \\\text{(by (\ref{sol.vander-det.xi+yj.short.c.detB}))}}%
}\cdot\underbrace{\det C}_{\substack{=p_{n-1}^{n}\left(  \prod_{k=0}%
^{n-1}\dbinom{n-1}{k}\right)  \\\text{(by
(\ref{sol.vander-det.xi+yj.short.c.detC}))}}}\cdot\underbrace{\det
D}_{\substack{=\prod_{1\leq i<j\leq n}\left(  y_{j}-y_{i}\right)  \\\text{(by
(\ref{sol.vander-det.xi+yj.short.c.detD}))}}}\\
&  =\left(  \prod_{1\leq i<j\leq n}\left(  x_{i}-x_{j}\right)  \right)  \cdot
p_{n-1}^{n}\left(  \prod_{k=0}^{n-1}\dbinom{n-1}{k}\right)  \cdot\left(
\prod_{1\leq i<j\leq n}\left(  y_{j}-y_{i}\right)  \right) \\
&  =p_{n-1}^{n}\left(  \prod_{k=0}^{n-1}\dbinom{n-1}{k}\right)  \left(
\prod_{1\leq i<j\leq n}\left(  x_{i}-x_{j}\right)  \right)  \left(
\prod_{1\leq i<j\leq n}\left(  y_{j}-y_{i}\right)  \right)  .
\end{align*}
This solves Exercise \ref{exe.vander-det.xi+yj} \textbf{(c)}.

\textbf{(d)} We define two $n\times n$-matrices $F$ and $G$:

\begin{itemize}
\item Define an $n\times n$-matrix $F$ by
\[
F=\left(  x_{i}^{j-1}\right)  _{1\leq i\leq n,\ 1\leq j\leq n}.
\]
Thus,
\begin{align*}
\det F  &  =\det\left(  \left(  x_{i}^{j-1}\right)  _{1\leq i\leq n,\ 1\leq
j\leq n}\right) \\
&  =\prod_{1\leq j<i\leq n}\left(  x_{i}-x_{j}\right)
\ \ \ \ \ \ \ \ \ \ \left(  \text{by Theorem \ref{thm.vander-det}
\textbf{(c)}}\right)  .
\end{align*}

\item Define an $n\times n$-matrix $G$ by%
\[
G=\left(  p_{i-1}y_{j}^{i-1}\right)  _{1\leq i\leq n,\ 1\leq j\leq n}.
\]
The definition of $\det G$ yields%
\begin{align*}
\det G  &  =\sum_{\sigma\in S_{n}}\left(  -1\right)  ^{\sigma}%
\underbrace{\left(  p_{1-1}y_{\sigma\left(  1\right)  }^{1-1}\right)  \left(
p_{2-1}y_{\sigma\left(  2\right)  }^{2-1}\right)  \cdots\left(  p_{n-1}%
y_{\sigma\left(  n\right)  }^{n-1}\right)  }_{\substack{=\prod_{k=1}%
^{n}\left(  p_{k-1}y_{\sigma\left(  k\right)  }^{k-1}\right)  \\=\left(
\prod_{k=1}^{n}p_{k-1}\right)  \left(  \prod_{k=1}^{n}y_{\sigma\left(
k\right)  }^{k-1}\right)  }}\\
&  =\sum_{\sigma\in S_{n}}\left(  -1\right)  ^{\sigma}\left(  \prod_{k=1}%
^{n}p_{k-1}\right)  \left(  \prod_{k=1}^{n}y_{\sigma\left(  k\right)  }%
^{k-1}\right) \\
&  =\left(  \prod_{k=1}^{n}p_{k-1}\right)  \cdot\sum_{\sigma\in S_{n}}\left(
-1\right)  ^{\sigma}\prod_{k=1}^{n}y_{\sigma\left(  k\right)  }^{k-1}.
\end{align*}
In light of%
\begin{align*}
&  \sum_{\sigma\in S_{n}}\left(  -1\right)  ^{\sigma}\underbrace{\prod
_{k=1}^{n}y_{\sigma\left(  k\right)  }^{k-1}}_{=y_{\sigma\left(  1\right)
}^{1-1}y_{\sigma\left(  2\right)  }^{2-1}\cdots y_{\sigma\left(  n\right)
}^{n-1}}\\
&  =\sum_{\sigma\in S_{n}}\left(  -1\right)  ^{\sigma}y_{\sigma\left(
1\right)  }^{1-1}y_{\sigma\left(  2\right)  }^{2-1}\cdots y_{\sigma\left(
n\right)  }^{n-1}=\det\left(  \left(  y_{j}^{i-1}\right)  _{1\leq i\leq
n,\ 1\leq j\leq n}\right) \\
&  \ \ \ \ \ \ \ \ \ \ \ \ \ \ \ \ \ \ \ \ \left(  \text{by the definition of
}\det\left(  \left(  y_{j}^{i-1}\right)  _{1\leq i\leq n,\ 1\leq j\leq
n}\right)  \right) \\
&  =\prod_{1\leq j<i\leq n}\left(  y_{i}-y_{j}\right) \\
&  \ \ \ \ \ \ \ \ \ \ \ \ \ \ \ \ \ \ \ \ \left(  \text{by Theorem
\ref{thm.vander-det} \textbf{(d)}, applied to }y_{k}\text{ instead of }%
x_{k}\right)  ,
\end{align*}
we can rewrite this as%
\[
\det G=\left(  \prod_{k=1}^{n}p_{k-1}\right)  \cdot\prod_{1\leq j<i\leq
n}\left(  y_{i}-y_{j}\right)  .
\]

\end{itemize}

We shall next prove that $\left(  P\left(  x_{i}y_{j}\right)  \right)  _{1\leq
i\leq n,\ 1\leq j\leq n}=FG$. Once this is done, we will be able to compute
$\det\left(  \left(  P\left(  x_{i}y_{j}\right)  \right)  _{1\leq i\leq
n,\ 1\leq j\leq n}\right)  $ by applying Theorem \ref{thm.det(AB)}.

Let $\left(  i,j\right)  \in\left\{  1,2,\ldots,n\right\}  ^{2}$. Then,
substituting $x_{i}y_{j}$ for $X$ on both sides of the equality $P\left(
X\right)  =\sum_{k=0}^{n-1}p_{k}X^{k}$, we obtain%
\[
P\left(  x_{i}y_{j}\right)  =\sum_{k=0}^{n-1}\underbrace{p_{k}\left(
x_{i}y_{j}\right)  ^{k}}_{\substack{=p_{k}x_{i}^{k}y_{j}^{k}\\=x_{i}^{k}%
p_{k}y_{j}^{k}}}=\sum_{k=0}^{n-1}x_{i}^{k}p_{k}y_{j}^{k}=\sum_{k=1}^{n}%
x_{i}^{k-1}p_{k-1}y_{j}^{k-1}%
\]
(here, we have substituted $k-1$ for $k$ in the sum).

Forget that we fixed $\left(  i,j\right)  $. We thus have shown that $P\left(
x_{i}y_{j}\right)  =\sum_{k=1}^{n}x_{i}^{k-1}p_{k-1}y_{j}^{k-1}$ for each
$\left(  i,j\right)  \in\left\{  1,2,\ldots,n\right\}  ^{2}$. In other words,%
\[
\left(  P\left(  x_{i}y_{j}\right)  \right)  _{1\leq i\leq n,\ 1\leq j\leq
n}=\left(  \sum_{k=1}^{n}x_{i}^{k-1}p_{k-1}y_{j}^{k-1}\right)  _{1\leq i\leq
n,\ 1\leq j\leq n}%
\]

On the other hand, the definition of the product of two matrices shows that%
\[
FG=\left(  \sum_{k=1}^{n}x_{i}^{k-1}p_{k-1}y_{j}^{k-1}\right)  _{1\leq i\leq
n,\ 1\leq j\leq n}%
\]
(since $F=\left(  x_{i}^{j-1}\right)  _{1\leq i\leq n,\ 1\leq j\leq n}$ and
$G=\left(  p_{i-1}y_{j}^{i-1}\right)  _{1\leq i\leq n,\ 1\leq j\leq n}$).
Comparing these two equalities, we obtain $\left(  P\left(  x_{i}y_{j}\right)
\right)  _{1\leq i\leq n,\ 1\leq j\leq n}=FG$. Therefore,%
\begin{align*}
&  \det\underbrace{\left(  \left(  P\left(  x_{i}y_{j}\right)  \right)
_{1\leq i\leq n,\ 1\leq j\leq n}\right)  }_{=FG}\\
&  =\det\left(  FG\right)  =\underbrace{\det F}_{=\prod_{1\leq j<i\leq
n}\left(  x_{i}-x_{j}\right)  }\cdot\underbrace{\det G}_{=\left(  \prod
_{k=1}^{n}p_{k-1}\right)  \cdot\prod_{1\leq j<i\leq n}\left(  y_{i}%
-y_{j}\right)  }\\
&  \ \ \ \ \ \ \ \ \ \ \left(  \text{by Theorem \ref{thm.det(AB)}, applied to
}F\text{ and }G\text{ instead of }A\text{ and }B\right) \\
&  =\left(  \prod_{1\leq j<i\leq n}\left(  x_{i}-x_{j}\right)  \right)
\cdot\left(  \prod_{k=1}^{n}p_{k-1}\right)  \cdot\prod_{1\leq j<i\leq
n}\left(  y_{i}-y_{j}\right) \\
&  =\underbrace{\left(  \prod_{k=1}^{n}p_{k-1}\right)  }_{\substack{=\prod
_{k=0}^{n-1}p_{k}\\\text{(here, we have}\\\text{substituted }k\\\text{for
}k-1\text{ in the product)}}}\underbrace{\left(  \prod_{1\leq j<i\leq
n}\left(  x_{i}-x_{j}\right)  \right)  \left(  \prod_{1\leq j<i\leq n}\left(
y_{i}-y_{j}\right)  \right)  }_{\substack{=\prod_{1\leq j<i\leq n}\left(
\left(  x_{i}-x_{j}\right)  \left(  y_{i}-y_{j}\right)  \right)
\\=\prod_{1\leq i<j\leq n}\left(  \left(  x_{j}-x_{i}\right)  \left(
y_{j}-y_{i}\right)  \right)  \\\text{(here, we have renamed the}\\\text{index
}\left(  j,i\right)  \text{ as }\left(  i,j\right)  \text{ in the product)}%
}}\\
&  =\left(  \prod_{k=0}^{n-1}p_{k}\right)  \prod_{1\leq i<j\leq n}%
\underbrace{\left(  \left(  x_{j}-x_{i}\right)  \left(  y_{j}-y_{i}\right)
\right)  }_{=\left(  x_{i}-x_{j}\right)  \left(  y_{i}-y_{j}\right)  }\\
&  =\left(  \prod_{k=0}^{n-1}p_{k}\right)  \underbrace{\prod_{1\leq i<j\leq
n}\left(  \left(  x_{i}-x_{j}\right)  \left(  y_{i}-y_{j}\right)  \right)
}_{=\left(  \prod_{1\leq i<j\leq n}\left(  x_{i}-x_{j}\right)  \right)
\left(  \prod_{1\leq i<j\leq n}\left(  y_{i}-y_{j}\right)  \right)  }\\
&  =\left(  \prod_{k=0}^{n-1}p_{k}\right)  \left(  \prod_{1\leq i<j\leq
n}\left(  x_{i}-x_{j}\right)  \right)  \left(  \prod_{1\leq i<j\leq n}\left(
y_{i}-y_{j}\right)  \right)  .
\end{align*}
This solves Exercise \ref{exe.vander-det.xi+yj} \textbf{(d)}.

\textbf{(a)} For any two objects $i$ and $j$, we define $\delta_{i,j}$ to be
the element $%
\begin{cases}
1, & \text{if }i=j;\\
0, & \text{if }i\neq j
\end{cases}
$ of $\mathbb{K}$.

Let $m\in\left\{  0,1,\ldots,n-1\right\}  $. (You are reading this right: We
are not requiring $m$ to belong to $\left\{  0,1,\ldots,n-2\right\}  $; the
purpose of this is to obtain a result that will bring us close to solving both
parts \textbf{(a)} and \textbf{(b)} simultaneously.)

Define an $n$-tuple $\left(  p_{0},p_{1},\ldots,p_{n-1}\right)  \in
\mathbb{K}^{n}$ by%
\[
\left(  p_{k}=\delta_{k,m}\ \ \ \ \ \ \ \ \ \ \text{for every }k\in\left\{
0,1,\ldots,n-1\right\}  \right)  .
\]
Let $P\left(  X\right)  \in\mathbb{K}\left[  X\right]  $ be the polynomial
$\sum_{k=0}^{n-1}p_{k}X^{k}$.

We have%
\begin{align*}
P\left(  X\right)   &  =\underbrace{\sum_{k=0}^{n-1}}_{=\sum_{k\in\left\{
0,1,\ldots,n-1\right\}  }}\underbrace{p_{k}}_{\substack{=\delta_{k,m}}%
}X^{k}=\sum_{k\in\left\{  0,1,\ldots,n-1\right\}  }\delta_{k,m}X^{k}\\
&  =\underbrace{\delta_{m,m}}_{\substack{=1\\\text{(since }m=m\text{)}}%
}X^{m}+\sum_{\substack{k\in\left\{  0,1,\ldots,n-1\right\}  ;\\k\neq
m}}\underbrace{\delta_{k,m}}_{\substack{=0\\\text{(since }k\neq m\text{)}%
}}X^{k}\\
&  \ \ \ \ \ \ \ \ \ \ \left(
\begin{array}
[c]{c}%
\text{here, we have split off the addend for }k=m\text{ from the sum}\\
\text{(since }m\in\left\{  0,1,\ldots,n-1\right\}  \text{)}%
\end{array}
\right) \\
&  =X^{m}+\underbrace{\sum_{\substack{k\in\left\{  0,1,\ldots,n-1\right\}
;\\k\neq m}}0X^{k}}_{=0}=X^{m}.
\end{align*}
Every $\left(  i,j\right)  \in\left\{  1,2,\ldots,n\right\}  ^{2}$ satisfies
$P\left(  x_{i}+y_{j}\right)  =\left(  x_{i}+y_{j}\right)  ^{m}$ (since
$P\left(  X\right)  =X^{m}$). In other words, $\left(  P\left(  x_{i}%
+y_{j}\right)  \right)  _{1\leq i\leq n,\ 1\leq j\leq n}=\left(  \left(
x_{i}+y_{j}\right)  ^{m}\right)  _{1\leq i\leq n,\ 1\leq j\leq n}$. Hence,%
\[
\left(  \left(  x_{i}+y_{j}\right)  ^{m}\right)  _{1\leq i\leq n,\ 1\leq j\leq
n}=\left(  P\left(  x_{i}+y_{j}\right)  \right)  _{1\leq i\leq n,\ 1\leq j\leq
n}.
\]
Taking determinants on both sides of this equality, we obtain%
\begin{align}
&  \det\left(  \left(  \left(  x_{i}+y_{j}\right)  ^{m}\right)  _{1\leq i\leq
n,\ 1\leq j\leq n}\right) \nonumber\\
&  =\det\left(  \left(  P\left(  x_{i}+y_{j}\right)  \right)  _{1\leq i\leq
n,\ 1\leq j\leq n}\right) \nonumber\\
&  =\underbrace{p_{n-1}^{n}}_{\substack{=\delta_{n-1,m}^{n}\\\text{(since
}p_{n-1}=\delta_{n-1,m}\\\text{(by the definition of }p_{n-1}\text{))}%
}}\left(  \prod_{k=0}^{n-1}\dbinom{n-1}{k}\right)  \left(  \prod_{1\leq
i<j\leq n}\left(  x_{i}-x_{j}\right)  \right)  \left(  \prod_{1\leq i<j\leq
n}\left(  y_{j}-y_{i}\right)  \right) \nonumber\\
&  \ \ \ \ \ \ \ \ \ \ \left(  \text{by Exercise \ref{exe.vander-det.xi+yj}
\textbf{(c)}}\right) \nonumber\\
&  =\delta_{n-1,m}^{n}\left(  \prod_{k=0}^{n-1}\dbinom{n-1}{k}\right)  \left(
\prod_{1\leq i<j\leq n}\left(  x_{i}-x_{j}\right)  \right)  \left(
\prod_{1\leq i<j\leq n}\left(  y_{j}-y_{i}\right)  \right)  .
\label{sol.vander-det.xi+yj.short.a.almost}%
\end{align}

Now, let us forget that we fixed $m$. We thus have proven
(\ref{sol.vander-det.xi+yj.short.a.almost}) for every $m\in\left\{
0,1,\ldots,n-1\right\}  $.

Now, let $m\in\left\{  0,1,\ldots,n-2\right\}  $. Thus, $m\neq n-1$, so that
$\delta_{n-1,m}=0$. Hence, $\delta_{n-1,m}^{n}=0^{n}=0$ (since $n$ is a
positive integer).

On the other hand, $m\in\left\{  0,1,\ldots,n-2\right\}  \subseteq\left\{
0,1,\ldots,n-1\right\}  $, and therefore
(\ref{sol.vander-det.xi+yj.short.a.almost}) holds. Hence,%
\begin{align*}
&  \det\left(  \left(  \left(  x_{i}+y_{j}\right)  ^{m}\right)  _{1\leq i\leq
n,\ 1\leq j\leq n}\right) \\
&  =\underbrace{\delta_{n-1,m}^{n}}_{=0}\left(  \prod_{k=0}^{n-1}\dbinom
{n-1}{k}\right)  \left(  \prod_{1\leq i<j\leq n}\left(  x_{i}-x_{j}\right)
\right)  \left(  \prod_{1\leq i<j\leq n}\left(  y_{j}-y_{i}\right)  \right) \\
&  =0.
\end{align*}
This solves Exercise \ref{exe.vander-det.xi+yj} \textbf{(a)}.

\textbf{(b)} For any two objects $i$ and $j$, we define $\delta_{i,j}$ as in
the solution to Exercise \ref{exe.vander-det.xi+yj} \textbf{(a)}.

We have $\delta_{n-1,n-1}=1$ and thus $\delta_{n-1,n-1}^{n}=1^{n}=1$.

In our solution of Exercise \ref{exe.vander-det.xi+yj} \textbf{(a)}, we have
proven the equality (\ref{sol.vander-det.xi+yj.short.a.almost}) for every
$m\in\left\{  0,1,\ldots,n-1\right\}  $. Thus, we can apply this equality to
$m=n-1$. As a result, we obtain
\begin{align*}
&  \det\left(  \left(  \left(  x_{i}+y_{j}\right)  ^{n-1}\right)  _{1\leq
i\leq n,\ 1\leq j\leq n}\right) \\
&  =\underbrace{\delta_{n-1,n-1}^{n}}_{=1}\left(  \prod_{k=0}^{n-1}%
\dbinom{n-1}{k}\right)  \left(  \prod_{1\leq i<j\leq n}\left(  x_{i}%
-x_{j}\right)  \right)  \left(  \prod_{1\leq i<j\leq n}\left(  y_{j}%
-y_{i}\right)  \right) \\
&  =\left(  \prod_{k=0}^{n-1}\dbinom{n-1}{k}\right)  \left(  \prod_{1\leq
i<j\leq n}\left(  x_{i}-x_{j}\right)  \right)  \left(  \prod_{1\leq i<j\leq
n}\left(  y_{j}-y_{i}\right)  \right)  .
\end{align*}
This solves Exercise \ref{exe.vander-det.xi+yj} \textbf{(b)}.
\end{proof}
\end{vershort}

\begin{verlong}
\begin{proof}
[Solution to Exercise \ref{exe.vander-det.xi+yj}.]\textbf{(c)} We extend the
$n$-tuple $\left(  p_{0},p_{1},\ldots,p_{n-1}\right)  $ to an infinite
sequence $\left(  p_{0},p_{1},p_{2},\ldots\right)  \in\mathbb{K}^{\infty}$ by
setting%
\begin{equation}
p_{\ell}=0\ \ \ \ \ \ \ \ \ \ \text{for every }\ell\in\left\{
n,n+1,n+2,\ldots\right\}  . \label{sol.vander-det.xi+yj.c.pk=0}%
\end{equation}

Next, we define three $n\times n$-matrices $B$, $C$ and $D$:

\begin{itemize}
\item Define an $n\times n$-matrix $B$ by
\[
B=\left(  x_{i}^{n-j}\right)  _{1\leq i\leq n,\ 1\leq j\leq n}.
\]
Then,%
\begin{align}
\det\underbrace{B}_{=\left(  x_{i}^{n-j}\right)  _{1\leq i\leq n,\ 1\leq j\leq
n}}  &  =\det\left(  \left(  x_{i}^{n-j}\right)  _{1\leq i\leq n,\ 1\leq j\leq
n}\right) \nonumber\\
&  =\prod_{1\leq i<j\leq n}\left(  x_{i}-x_{j}\right)
\label{sol.vander-det.xi+yj.c.detB}%
\end{align}
(by Theorem \ref{thm.vander-det} \textbf{(a)}).

\item For every $\left(  i,j\right)  \in\left\{  1,2,\ldots,n\right\}  ^{2}$,
the element $\dbinom{n-1-i+j}{j-1}p_{n-1-i+j}$ of $\mathbb{K}$ is
well-defined\footnote{\textit{Proof.} Let $\left(  i,j\right)  \in\left\{
1,2,\ldots,n\right\}  ^{2}$. Thus, $i\in\left\{  1,2,\ldots,n\right\}  $ and
$j\in\left\{  1,2,\ldots,n\right\}  $. Hence, $1\leq i\leq n$ (since
$i\in\left\{  1,2,\ldots,n\right\}  $) and $1\leq j\leq n$ (since
$j\in\left\{  1,2,\ldots,n\right\}  $). From $1\leq j$, we obtain $j-1\geq0$.
Hence, the binomial coefficient $\dbinom{n-1-i+j}{j-1}\in\mathbb{Z}$ is
well-defined.
\par
Moreover, $n-1-\underbrace{i}_{\leq n}+\underbrace{j}_{\geq1}\geq n-1-n+1=0$,
so that $n-1-i+j\in\mathbb{N}$. Thus, $p_{n-1-i+j}\in\mathbb{K}$ is
well-defined (since we have an infinite sequence $\left(  p_{0},p_{1}%
,p_{2},\ldots\right)  \in\mathbb{K}^{\infty}$). Hence, the element
$\dbinom{n-1-i+j}{j-1}p_{n-1-i+j}$ of $\mathbb{K}$ is well-defined. Qed.}.
Thus, for every $\left(  i,j\right)  \in\left\{  1,2,\ldots,n\right\}  ^{2}$,
we can define an element $c_{i,j}\in\mathbb{K}$ by $c_{i,j}=\dbinom
{n-1-i+j}{j-1}p_{n-1-i+j}$. Consider these elements $c_{i,j}$. Define an
$n\times n$-matrix $C$ by%
\[
C=\left(  c_{i,j}\right)  _{1\leq i\leq n,\ 1\leq j\leq n}.
\]
Consider this matrix $C$. We have $c_{i,j}=0$ for every $\left(  i,j\right)
\in\left\{  1,2,\ldots,n\right\}  ^{2}$ satisfying $i<j$%
\ \ \ \ \footnote{\textit{Proof.} Let $\left(  i,j\right)  \in\left\{
1,2,\ldots,n\right\}  ^{2}$ be such that $i<j$. From $i<j$, we obtain $i\leq
j-1$ (since $i$ and $j$ are integers). Thus, $n-1-\underbrace{i}_{\leq
j-1}+j\geq n-1-\left(  j-1\right)  +j=n$. In other words, $n-1-i+j\in\left\{
n,n+1,n+2,\ldots\right\}  $. Hence, $p_{n-1-i+j}=0$ (by
(\ref{sol.vander-det.xi+yj.c.pk=0}), applied to $\ell=n-1-i+j$). Hence,
$c_{i,j}=\dbinom{n-1-i+j}{j-1}\underbrace{p_{n-1-i+j}}_{=0}=0$, qed.}. Hence,
Exercise \ref{exe.ps4.3} (applied to $C$ and $\left(  c_{i,j}\right)  _{1\leq
i\leq n,\ 1\leq j\leq n}$ instead of $A$ and $\left(  a_{i,j}\right)  _{1\leq
i\leq n,\ 1\leq j\leq n}$) shows that%
\begin{align}
\det C  &  =c_{1,1}c_{2,2}\cdots c_{n,n}=\prod_{i=1}^{n}\underbrace{c_{i,i}%
}_{\substack{=\dbinom{n-1-i+i}{i-1}p_{n-1-i+i}\\\text{(by the definition of
}c_{i,i}\text{)}}}\nonumber\\
&  =\prod_{i=1}^{n}\left(  \underbrace{\dbinom{n-1-i+i}{i-1}}_{=\dbinom
{n-1}{i-1}}\underbrace{p_{n-1-i+i}}_{\substack{=p_{n-1}\\\text{(since
}n-1-i+i=n-1\text{)}}}\right) \nonumber\\
&  =\prod_{i=1}^{n}\left(  \dbinom{n-1}{i-1}p_{n-1}\right)  =\left(
\prod_{i=1}^{n}\dbinom{n-1}{i-1}\right)  p_{n-1}^{n}\nonumber\\
&  =p_{n-1}^{n}\left(  \prod_{i=1}^{n}\dbinom{n-1}{i-1}\right)  =p_{n-1}%
^{n}\left(  \prod_{k=0}^{n-1}\dbinom{n-1}{k}\right)
\label{sol.vander-det.xi+yj.c.detC}%
\end{align}
(here, we have substituted $k$ for $i-1$ in the product).

\item Define an $n\times n$-matrix $D$ by
\[
D=\left(  y_{j}^{i-1}\right)  _{1\leq i\leq n,\ 1\leq j\leq n}.
\]
Then,%
\begin{align}
\det\underbrace{D}_{=\left(  y_{j}^{i-1}\right)  _{1\leq i\leq n,\ 1\leq j\leq
n}}  &  =\det\left(  \left(  y_{j}^{i-1}\right)  _{1\leq i\leq n,\ 1\leq j\leq
n}\right)  =\prod_{1\leq j<i\leq n}\left(  y_{i}-y_{j}\right) \nonumber\\
&  \ \ \ \ \ \ \ \ \ \ \left(  \text{by Theorem \ref{thm.vander-det}
\textbf{(d)}, applied to }y_{k}\text{ instead of }x_{k}\right) \nonumber\\
&  =\prod_{1\leq i<j\leq n}\left(  y_{j}-y_{i}\right)
\label{sol.vander-det.xi+yj.c.detD}%
\end{align}
(here, we have renamed the index $\left(  j,i\right)  $ as $\left(
i,j\right)  $ in the product).
\end{itemize}

Our goal is to prove that $\left(  P\left(  x_{i}+y_{j}\right)  \right)
_{1\leq i\leq n,\ 1\leq j\leq n}=BCD$. Once this is done, we will be able to
compute $\det\left(  \left(  P\left(  x_{i}+y_{j}\right)  \right)  _{1\leq
i\leq n,\ 1\leq j\leq n}\right)  $ by applying Theorem \ref{thm.det(AB)} twice.

We have $C=\left(  c_{i,j}\right)  _{1\leq i\leq n,\ 1\leq j\leq n}$ and
$D=\left(  y_{j}^{i-1}\right)  _{1\leq i\leq n,\ 1\leq j\leq n}$. Hence, the
definition of the product of two matrices shows that%
\begin{equation}
CD=\left(  \sum_{k=1}^{n}c_{i,k}y_{j}^{k-1}\right)  _{1\leq i\leq n,\ 1\leq
j\leq n}. \label{sol.vander-det.xi+yj.c.CD.1}%
\end{equation}
But every $\left(  i,j\right)  \in\left\{  1,2,\ldots,n\right\}  ^{2}$
satisfies%
\begin{align*}
\sum_{k=1}^{n}c_{i,k}y_{j}^{k-1}  &  =\sum_{k=0}^{n-1}\underbrace{c_{i,k+1}%
}_{\substack{=\dbinom{n-1-i+\left(  k+1\right)  }{\left(  k+1\right)
-1}p_{n-1-i+\left(  k+1\right)  }\\\text{(by the definition of }%
c_{i,k+1}\text{)}}}y_{j}^{\left(  k+1\right)  -1}\\
&  \ \ \ \ \ \ \ \ \ \ \left(  \text{here, we have substituted }k+1\text{ for
}k\text{ in the sum}\right) \\
&  =\sum_{k=0}^{n-1}\underbrace{\dbinom{n-1-i+\left(  k+1\right)  }{\left(
k+1\right)  -1}}_{\substack{=\dbinom{n-i+k}{k}\\\text{(since }n-1-i+\left(
k+1\right)  =n-i+k\\\text{and }\left(  k+1\right)  -1=k\text{)}}%
}\underbrace{p_{n-1-i+\left(  k+1\right)  }}_{\substack{=p_{n-i+k}%
\\\text{(since }n-1-i+\left(  k+1\right)  =n-i+k\text{)}}}\underbrace{y_{j}%
^{\left(  k+1\right)  -1}}_{\substack{=y_{j}^{k}\\\text{(since }\left(
k+1\right)  -1=k\text{)}}}\\
&  =\sum_{k=0}^{n-1}\dbinom{n-i+k}{k}p_{n-i+k}y_{j}^{k}=\sum_{\ell=0}%
^{n-1}\dbinom{n-i+\ell}{\ell}p_{n-i+\ell}y_{j}^{\ell}%
\end{align*}
(here, we renamed the summation index $k$ as $\ell$). Hence,
(\ref{sol.vander-det.xi+yj.c.CD.1}) becomes%
\begin{align*}
CD  &  =\left(  \underbrace{\sum_{k=1}^{n}c_{i,k}y_{j}^{k-1}}_{=\sum_{\ell
=0}^{n-1}\dbinom{n-i+\ell}{\ell}p_{n-i+\ell}y_{j}^{\ell}}\right)  _{1\leq
i\leq n,\ 1\leq j\leq n}\\
&  =\left(  \sum_{\ell=0}^{n-1}\dbinom{n-i+\ell}{\ell}p_{n-i+\ell}y_{j}^{\ell
}\right)  _{1\leq i\leq n,\ 1\leq j\leq n}.
\end{align*}

Now, $B=\left(  x_{i}^{n-j}\right)  _{1\leq i\leq n,\ 1\leq j\leq n}$ and
$CD=\left(  \sum_{\ell=0}^{n-1}\dbinom{n-i+\ell}{\ell}p_{n-i+\ell}y_{j}^{\ell
}\right)  _{1\leq i\leq n,\ 1\leq j\leq n}$. Hence, the definition of the
product of two matrices shows that%
\begin{equation}
B\left(  CD\right)  =\left(  \sum_{k=1}^{n}x_{i}^{n-k}\left(  \sum_{\ell
=0}^{n-1}\dbinom{n-k+\ell}{\ell}p_{n-k+\ell}y_{j}^{\ell}\right)  \right)
_{1\leq i\leq n,\ 1\leq j\leq n}. \label{sol.vander-det.xi+yj.c.BCD.1}%
\end{equation}

Now, let $\left(  i,j\right)  \in\left\{  1,2,\ldots,n\right\}  ^{2}$. Then,%
\begin{align}
&  \sum_{k=1}^{n}x_{i}^{n-k}\left(  \sum_{\ell=0}^{n-1}\dbinom{n-k+\ell}{\ell
}p_{n-k+\ell}y_{j}^{\ell}\right) \nonumber\\
&  =\sum_{k=0}^{n-1}x_{i}^{k}\left(  \sum_{\ell=0}^{n-1}\dbinom{k+\ell}{\ell
}p_{k+\ell}y_{j}^{\ell}\right) \nonumber\\
&  \ \ \ \ \ \ \ \ \ \ \left(  \text{here, we have substituted }k\text{ for
}n-k\text{ in the first sum}\right) \nonumber\\
&  =\sum_{k=0}^{n-1}\sum_{\ell=0}^{n-1}x_{i}^{k}\dbinom{k+\ell}{\ell}%
p_{k+\ell}y_{j}^{\ell}\nonumber\\
&  =\sum_{k=0}^{n-1}\sum_{\ell=k}^{n-1+k}x_{i}^{k}\underbrace{\dbinom
{k+\left(  \ell-k\right)  }{\ell-k}}_{\substack{=\dbinom{\ell}{\ell
-k}\\\text{(since }k+\left(  \ell-k\right)  =\ell\text{)}}%
}\underbrace{p_{k+\left(  \ell-k\right)  }}_{\substack{=p_{\ell}\\\text{(since
}k+\left(  \ell-k\right)  =\ell\text{)}}}y_{j}^{\ell-k}\nonumber\\
&  \ \ \ \ \ \ \ \ \ \ \left(  \text{here, we have substituted }\ell-k\text{
for }\ell\text{ in the second sum}\right) \nonumber\\
&  =\sum_{k=0}^{n-1}\underbrace{\sum_{\ell=k}^{n-1+k}x_{i}^{k}\dbinom{\ell
}{\ell-k}p_{\ell}y_{j}^{\ell-k}}_{\substack{=\sum_{\ell=k}^{n-1}x_{i}%
^{k}\dbinom{\ell}{\ell-k}p_{\ell}y_{j}^{\ell-k}+\sum_{\ell=n}^{n-1+k}x_{i}%
^{k}\dbinom{\ell}{\ell-k}p_{\ell}y_{j}^{\ell-k}\\\text{(since }k\leq n\leq
n-1+k+1\text{ (since }n\leq n+k=n-1+k+1\text{))}}}\nonumber\\
&  =\sum_{k=0}^{n-1}\left(  \sum_{\ell=k}^{n-1}x_{i}^{k}\dbinom{\ell}{\ell
-k}p_{\ell}y_{j}^{\ell-k}+\sum_{\ell=n}^{n-1+k}x_{i}^{k}\dbinom{\ell}{\ell
-k}\underbrace{p_{\ell}}_{\substack{=0\\\text{(by
(\ref{sol.vander-det.xi+yj.c.pk=0})}\\\text{(since }\ell\in\left\{
n,n+1,\ldots,n-1+k\right\}  \\\subseteq\left\{  n,n+1,n+2,\ldots\right\}
\text{))}}}y_{j}^{\ell-k}\right) \nonumber\\
&  =\sum_{k=0}^{n-1}\left(  \sum_{\ell=k}^{n-1}x_{i}^{k}\dbinom{\ell}{\ell
-k}p_{\ell}y_{j}^{\ell-k}+\underbrace{\sum_{\ell=n}^{n-1+k}x_{i}^{k}%
\dbinom{\ell}{\ell-k}0y_{j}^{\ell-k}}_{=0}\right) \nonumber\\
&  =\sum_{k=0}^{n-1}\sum_{\ell=k}^{n-1}x_{i}^{k}\dbinom{\ell}{\ell-k}p_{\ell
}y_{j}^{\ell-k}. \label{sol.vander-det.xi+yj.c.BCD.2}%
\end{align}
On the other hand, $P\left(  X\right)  =\sum_{k=0}^{n-1}p_{k}X^{k}=\sum
_{\ell=0}^{n-1}p_{\ell}X^{\ell}$ (here, we renamed the summation index $k$ as
$\ell$). Hence,%
\begin{align}
P\left(  x_{i}+y_{j}\right)   &  =\sum_{\ell=0}^{n-1}p_{\ell}%
\underbrace{\left(  x_{i}+y_{j}\right)  ^{\ell}}_{\substack{=\sum_{k=0}^{\ell
}\dbinom{\ell}{k}x_{i}^{k}y_{j}^{\ell-k}\\\text{(by the binomial formula)}%
}}=\sum_{\ell=0}^{n-1}p_{\ell}\sum_{k=0}^{\ell}\dbinom{\ell}{k}x_{i}^{k}%
y_{j}^{\ell-k}\nonumber\\
&  =\underbrace{\sum_{\ell=0}^{n-1}\sum_{k=0}^{\ell}}_{=\sum
_{\substack{\left(  \ell,k\right)  \in\left\{  0,1,\ldots,n-1\right\}
^{2};\\k\leq\ell}}=\sum_{k=0}^{n-1}\sum_{\ell=k}^{n-1}}p_{\ell}\dbinom{\ell
}{k}x_{i}^{k}y_{j}^{\ell-k}\nonumber\\
&  =\sum_{k=0}^{n-1}\sum_{\ell=k}^{n-1}p_{\ell}\underbrace{\dbinom{\ell}{k}%
}_{\substack{=\dbinom{\ell}{\ell-k}\\\text{(by (\ref{eq.binom.symm}) (applied
to }\ell\text{ and }k\\\text{instead of }m\text{ and }n\text{))}}}x_{i}%
^{k}y_{j}^{\ell-k}\nonumber\\
&  =\sum_{k=0}^{n-1}\sum_{\ell=k}^{n-1}\underbrace{p_{\ell}\dbinom{\ell}%
{\ell-k}x_{i}^{k}}_{=x_{i}^{k}\dbinom{\ell}{\ell-k}p_{\ell}}y_{j}^{\ell
-k}=\sum_{k=0}^{n-1}\sum_{\ell=k}^{n-1}x_{i}^{k}\dbinom{\ell}{\ell-k}p_{\ell
}y_{j}^{\ell-k}\nonumber\\
&  =\sum_{k=1}^{n}x_{i}^{n-k}\left(  \sum_{\ell=0}^{n-1}\dbinom{n-k+\ell}%
{\ell}p_{n-k+\ell}y_{j}^{\ell}\right)  \ \ \ \ \ \ \ \ \ \ \left(  \text{by
(\ref{sol.vander-det.xi+yj.c.BCD.2})}\right)  .
\label{sol.vander-det.xi+yj.c.BCD.3}%
\end{align}

Now, let us forget that we fixed $\left(  i,j\right)  $. We thus have proven
(\ref{sol.vander-det.xi+yj.c.BCD.3}) for every $\left(  i,j\right)
\in\left\{  1,2,\ldots,n\right\}  ^{2}$. Hence,%
\begin{align*}
&  \left(  \underbrace{P\left(  x_{i}+y_{j}\right)  }_{\substack{=\sum
_{k=1}^{n}x_{i}^{n-k}\left(  \sum_{\ell=0}^{n-1}\dbinom{n-k+\ell}{\ell
}p_{n-k+\ell}y_{j}^{\ell}\right)  \\\text{(by
(\ref{sol.vander-det.xi+yj.c.BCD.3}))}}}\right)  _{1\leq i\leq n,\ 1\leq j\leq
n}\\
&  =\left(  \sum_{k=1}^{n}x_{i}^{n-k}\left(  \sum_{\ell=0}^{n-1}%
\dbinom{n-k+\ell}{\ell}p_{n-k+\ell}y_{j}^{\ell}\right)  \right)  _{1\leq i\leq
n,\ 1\leq j\leq n}=B\left(  CD\right)  \ \ \ \ \ \ \ \ \ \ \left(  \text{by
(\ref{sol.vander-det.xi+yj.c.BCD.1})}\right) \\
&  =BCD.
\end{align*}
Hence,%
\begin{align*}
&  \det\underbrace{\left(  \left(  P\left(  x_{i}+y_{j}\right)  \right)
_{1\leq i\leq n,\ 1\leq j\leq n}\right)  }_{=BCD}\\
&  =\det\left(  BCD\right)  =\det B\cdot\underbrace{\det\left(  CD\right)
}_{\substack{=\det C\cdot\det D\\\text{(by Theorem \ref{thm.det(AB)}, applied
to }C\text{ and }D\\\text{instead of }A\text{ and }B\text{)}}}\\
&  \ \ \ \ \ \ \ \ \ \ \left(  \text{by Theorem \ref{thm.det(AB)}, applied to
}B\text{ and }CD\text{ instead of }A\text{ and }B\right) \\
&  =\underbrace{\det B}_{\substack{=\prod_{1\leq i<j\leq n}\left(  x_{i}%
-x_{j}\right)  \\\text{(by (\ref{sol.vander-det.xi+yj.c.detB}))}}%
}\cdot\underbrace{\det C}_{\substack{=p_{n-1}^{n}\left(  \prod_{k=0}%
^{n-1}\dbinom{n-1}{k}\right)  \\\text{(by (\ref{sol.vander-det.xi+yj.c.detC}%
))}}}\cdot\underbrace{\det D}_{\substack{=\prod_{1\leq i<j\leq n}\left(
y_{j}-y_{i}\right)  \\\text{(by (\ref{sol.vander-det.xi+yj.c.detD}))}}}\\
&  =\left(  \prod_{1\leq i<j\leq n}\left(  x_{i}-x_{j}\right)  \right)  \cdot
p_{n-1}^{n}\left(  \prod_{k=0}^{n-1}\dbinom{n-1}{k}\right)  \cdot\left(
\prod_{1\leq i<j\leq n}\left(  y_{j}-y_{i}\right)  \right) \\
&  =p_{n-1}^{n}\left(  \prod_{k=0}^{n-1}\dbinom{n-1}{k}\right)  \left(
\prod_{1\leq i<j\leq n}\left(  x_{i}-x_{j}\right)  \right)  \left(
\prod_{1\leq i<j\leq n}\left(  y_{j}-y_{i}\right)  \right)  .
\end{align*}
This solves Exercise \ref{exe.vander-det.xi+yj} \textbf{(c)}.

\textbf{(d)} We define two $n\times n$-matrices $F$ and $G$:

\begin{itemize}
\item Define an $n\times n$-matrix $F$ by
\[
F=\left(  x_{i}^{j-1}\right)  _{1\leq i\leq n,\ 1\leq j\leq n}.
\]
Then,
\begin{align*}
\det\underbrace{F}_{=\left(  x_{i}^{j-1}\right)  _{1\leq i\leq n,\ 1\leq j\leq
n}}  &  =\det\left(  \left(  x_{i}^{j-1}\right)  _{1\leq i\leq n,\ 1\leq j\leq
n}\right) \\
&  =\prod_{1\leq j<i\leq n}\left(  x_{i}-x_{j}\right)
\ \ \ \ \ \ \ \ \ \ \left(  \text{by Theorem \ref{thm.vander-det}
\textbf{(c)}}\right)  .
\end{align*}

\item Define an $n\times n$-matrix $G$ by%
\[
G=\left(  p_{i-1}y_{j}^{i-1}\right)  _{1\leq i\leq n,\ 1\leq j\leq n}.
\]
The definition of $\det G$ yields%
\begin{align*}
\det G  &  =\sum_{\sigma\in S_{n}}\left(  -1\right)  ^{\sigma}%
\underbrace{\left(  p_{1-1}y_{\sigma\left(  1\right)  }^{1-1}\right)  \left(
p_{2-1}y_{\sigma\left(  2\right)  }^{2-1}\right)  \cdots\left(  p_{n-1}%
y_{\sigma\left(  n\right)  }^{n-1}\right)  }_{\substack{=\prod_{k=1}%
^{n}\left(  p_{k-1}y_{\sigma\left(  k\right)  }^{k-1}\right)  \\=\left(
\prod_{k=1}^{n}p_{k-1}\right)  \left(  \prod_{k=1}^{n}y_{\sigma\left(
k\right)  }^{k-1}\right)  }}\\
&  =\sum_{\sigma\in S_{n}}\left(  -1\right)  ^{\sigma}\left(  \prod_{k=1}%
^{n}p_{k-1}\right)  \left(  \prod_{k=1}^{n}y_{\sigma\left(  k\right)  }%
^{k-1}\right) \\
&  =\left(  \prod_{k=1}^{n}p_{k-1}\right)  \cdot\sum_{\sigma\in S_{n}}\left(
-1\right)  ^{\sigma}\prod_{k=1}^{n}y_{\sigma\left(  k\right)  }^{k-1}.
\end{align*}
In light of%
\begin{align*}
&  \sum_{\sigma\in S_{n}}\left(  -1\right)  ^{\sigma}\underbrace{\prod
_{k=1}^{n}y_{\sigma\left(  k\right)  }^{k-1}}_{=y_{\sigma\left(  1\right)
}^{1-1}y_{\sigma\left(  2\right)  }^{2-1}\cdots y_{\sigma\left(  n\right)
}^{n-1}}\\
&  =\sum_{\sigma\in S_{n}}\left(  -1\right)  ^{\sigma}y_{\sigma\left(
1\right)  }^{1-1}y_{\sigma\left(  2\right)  }^{2-1}\cdots y_{\sigma\left(
n\right)  }^{n-1}=\det\left(  \left(  y_{j}^{i-1}\right)  _{1\leq i\leq
n,\ 1\leq j\leq n}\right) \\
&  \ \ \ \ \ \ \ \ \ \ \ \ \ \ \ \ \ \ \ \ \left(
\begin{array}
[c]{c}%
\text{since the definition of }\det\left(  \left(  y_{j}^{i-1}\right)  _{1\leq
i\leq n,\ 1\leq j\leq n}\right) \\
\text{yields }\det\left(  \left(  y_{j}^{i-1}\right)  _{1\leq i\leq n,\ 1\leq
j\leq n}\right)  =\sum_{\sigma\in S_{n}}\left(  -1\right)  ^{\sigma}%
y_{\sigma\left(  1\right)  }^{1-1}y_{\sigma\left(  2\right)  }^{2-1}\cdots
y_{\sigma\left(  n\right)  }^{n-1}%
\end{array}
\right) \\
&  =\prod_{1\leq j<i\leq n}\left(  y_{i}-y_{j}\right) \\
&  \ \ \ \ \ \ \ \ \ \ \ \ \ \ \ \ \ \ \ \ \left(  \text{by Theorem
\ref{thm.vander-det} \textbf{(d)}, applied to }y_{k}\text{ instead of }%
x_{k}\right)  ,
\end{align*}
we can rewrite this as%
\[
\det G=\left(  \prod_{k=1}^{n}p_{k-1}\right)  \cdot\prod_{1\leq j<i\leq
n}\left(  y_{i}-y_{j}\right)  .
\]

\end{itemize}

We shall next prove that $\left(  P\left(  x_{i}y_{j}\right)  \right)  _{1\leq
i\leq n,\ 1\leq j\leq n}=FG$. Once this is done, we will be able to compute
$\det\left(  \left(  P\left(  x_{i}y_{j}\right)  \right)  _{1\leq i\leq
n,\ 1\leq j\leq n}\right)  $ by applying Theorem \ref{thm.det(AB)}.

Let $\left(  i,j\right)  \in\left\{  1,2,\ldots,n\right\}  ^{2}$. Then,
substituting $x_{i}y_{j}$ for $X$ on both sides of the equality $P\left(
X\right)  =\sum_{k=0}^{n-1}p_{k}X^{k}$, we obtain%
\[
P\left(  x_{i}y_{j}\right)  =\sum_{k=0}^{n-1}\underbrace{p_{k}\left(
x_{i}y_{j}\right)  ^{k}}_{\substack{=p_{k}x_{i}^{k}y_{j}^{k}\\=x_{i}^{k}%
p_{k}y_{j}^{k}}}=\sum_{k=0}^{n-1}x_{i}^{k}p_{k}y_{j}^{k}=\sum_{k=1}^{n}%
x_{i}^{k-1}p_{k-1}y_{j}^{k-1}%
\]
(here, we have substituted $k-1$ for $k$ in the sum).

Forget that we fixed $\left(  i,j\right)  $. We thus have shown that $P\left(
x_{i}y_{j}\right)  =\sum_{k=1}^{n}x_{i}^{k-1}p_{k-1}y_{j}^{k-1}$ for each
$\left(  i,j\right)  \in\left\{  1,2,\ldots,n\right\}  ^{2}$. In other words,%
\[
\left(  P\left(  x_{i}y_{j}\right)  \right)  _{1\leq i\leq n,\ 1\leq j\leq
n}=\left(  \sum_{k=1}^{n}x_{i}^{k-1}p_{k-1}y_{j}^{k-1}\right)  _{1\leq i\leq
n,\ 1\leq j\leq n}%
\]

On the other hand, the definition of the product of two matrices shows that%
\[
FG=\left(  \sum_{k=1}^{n}x_{i}^{k-1}p_{k-1}y_{j}^{k-1}\right)  _{1\leq i\leq
n,\ 1\leq j\leq n}%
\]
(since $F=\left(  x_{i}^{j-1}\right)  _{1\leq i\leq n,\ 1\leq j\leq n}$ and
$G=\left(  p_{i-1}y_{j}^{i-1}\right)  _{1\leq i\leq n,\ 1\leq j\leq n}$).
Comparing these two equalities, we obtain $\left(  P\left(  x_{i}y_{j}\right)
\right)  _{1\leq i\leq n,\ 1\leq j\leq n}=FG$. Therefore,%
\begin{align*}
&  \det\underbrace{\left(  \left(  P\left(  x_{i}y_{j}\right)  \right)
_{1\leq i\leq n,\ 1\leq j\leq n}\right)  }_{=FG}\\
&  =\det\left(  FG\right)  =\underbrace{\det F}_{=\prod_{1\leq j<i\leq
n}\left(  x_{i}-x_{j}\right)  }\cdot\underbrace{\det G}_{=\left(  \prod
_{k=1}^{n}p_{k-1}\right)  \cdot\prod_{1\leq j<i\leq n}\left(  y_{i}%
-y_{j}\right)  }\\
&  \ \ \ \ \ \ \ \ \ \ \left(  \text{by Theorem \ref{thm.det(AB)}, applied to
}F\text{ and }G\text{ instead of }A\text{ and }B\right) \\
&  =\left(  \prod_{1\leq j<i\leq n}\left(  x_{i}-x_{j}\right)  \right)
\cdot\left(  \prod_{k=1}^{n}p_{k-1}\right)  \cdot\prod_{1\leq j<i\leq
n}\left(  y_{i}-y_{j}\right) \\
&  =\underbrace{\left(  \prod_{k=1}^{n}p_{k-1}\right)  }_{\substack{=\prod
_{k=0}^{n-1}p_{k}\\\text{(here, we have}\\\text{substituted }k\\\text{for
}k-1\text{ in the product)}}}\underbrace{\left(  \prod_{1\leq j<i\leq
n}\left(  x_{i}-x_{j}\right)  \right)  \left(  \prod_{1\leq j<i\leq n}\left(
y_{i}-y_{j}\right)  \right)  }_{\substack{=\prod_{1\leq j<i\leq n}\left(
\left(  x_{i}-x_{j}\right)  \left(  y_{i}-y_{j}\right)  \right)
\\=\prod_{1\leq i<j\leq n}\left(  \left(  x_{j}-x_{i}\right)  \left(
y_{j}-y_{i}\right)  \right)  \\\text{(here, we have renamed the}\\\text{index
}\left(  j,i\right)  \text{ as }\left(  i,j\right)  \text{ in the product)}%
}}\\
&  =\left(  \prod_{k=0}^{n-1}p_{k}\right)  \prod_{1\leq i<j\leq n}%
\underbrace{\left(  \left(  x_{j}-x_{i}\right)  \left(  y_{j}-y_{i}\right)
\right)  }_{=\left(  x_{i}-x_{j}\right)  \left(  y_{i}-y_{j}\right)  }\\
&  =\left(  \prod_{k=0}^{n-1}p_{k}\right)  \underbrace{\prod_{1\leq i<j\leq
n}\left(  \left(  x_{i}-x_{j}\right)  \left(  y_{i}-y_{j}\right)  \right)
}_{=\left(  \prod_{1\leq i<j\leq n}\left(  x_{i}-x_{j}\right)  \right)
\left(  \prod_{1\leq i<j\leq n}\left(  y_{i}-y_{j}\right)  \right)  }\\
&  =\left(  \prod_{k=0}^{n-1}p_{k}\right)  \left(  \prod_{1\leq i<j\leq
n}\left(  x_{i}-x_{j}\right)  \right)  \left(  \prod_{1\leq i<j\leq n}\left(
y_{i}-y_{j}\right)  \right)  .
\end{align*}
This solves Exercise \ref{exe.vander-det.xi+yj} \textbf{(d)}.

\textbf{(a)} For any two objects $i$ and $j$, we define $\delta_{i,j}$ to be
the element $%
\begin{cases}
1, & \text{if }i=j;\\
0, & \text{if }i\neq j
\end{cases}
$ of $\mathbb{K}$.

Let $m\in\left\{  0,1,\ldots,n-1\right\}  $. (You are reading this right: We
are not requiring $m$ to belong to $\left\{  0,1,\ldots,n-2\right\}  $; the
purpose of this is to obtain a result that will bring us close to solving both
parts \textbf{(a)} and \textbf{(b)} simultaneously.)

Define an $n$-tuple $\left(  p_{0},p_{1},\ldots,p_{n-1}\right)  \in
\mathbb{K}^{n}$ by%
\[
\left(  p_{k}=\delta_{k,m}\ \ \ \ \ \ \ \ \ \ \text{for every }k\in\left\{
0,1,\ldots,n-1\right\}  \right)  .
\]
Let $P\left(  X\right)  \in\mathbb{K}\left[  X\right]  $ be the polynomial
$\sum_{k=0}^{n-1}p_{k}X^{k}$.

We have%
\begin{align*}
P\left(  X\right)   &  =\underbrace{\sum_{k=0}^{n-1}}_{=\sum_{k\in\left\{
0,1,\ldots,n-1\right\}  }}\underbrace{p_{k}}_{\substack{=\delta_{k,m}%
\\\text{(by the definition}\\\text{of }p_{k}\text{)}}}X^{k}=\sum_{k\in\left\{
0,1,\ldots,n-1\right\}  }\delta_{k,m}X^{k}\\
&  =\underbrace{\delta_{m,m}}_{\substack{=1\\\text{(since }m=m\text{)}}%
}X^{m}+\sum_{\substack{k\in\left\{  0,1,\ldots,n-1\right\}  ;\\k\neq
m}}\underbrace{\delta_{k,m}}_{\substack{=0\\\text{(since }k\neq m\text{)}%
}}X^{k}\\
&  \ \ \ \ \ \ \ \ \ \ \left(
\begin{array}
[c]{c}%
\text{here, we have split off the addend for }k=m\text{ from the sum}\\
\text{(since }m\in\left\{  0,1,\ldots,n-1\right\}  \text{)}%
\end{array}
\right) \\
&  =X^{m}+\underbrace{\sum_{\substack{k\in\left\{  0,1,\ldots,n-1\right\}
;\\k\neq m}}0X^{k}}_{=0}=X^{m}.
\end{align*}
Every $\left(  i,j\right)  \in\left\{  1,2,\ldots,n\right\}  ^{2}$ satisfies
$P\left(  x_{i}+y_{j}\right)  =\left(  x_{i}+y_{j}\right)  ^{m}$ (since
$P\left(  X\right)  =X^{m}$). In other words, $\left(  P\left(  x_{i}%
+y_{j}\right)  \right)  _{1\leq i\leq n,\ 1\leq j\leq n}=\left(  \left(
x_{i}+y_{j}\right)  ^{m}\right)  _{1\leq i\leq n,\ 1\leq j\leq n}$. Hence,%
\[
\left(  \left(  x_{i}+y_{j}\right)  ^{m}\right)  _{1\leq i\leq n,\ 1\leq j\leq
n}=\left(  P\left(  x_{i}+y_{j}\right)  \right)  _{1\leq i\leq n,\ 1\leq j\leq
n}.
\]
Taking determinants on both sides of this equality, we obtain%
\begin{align}
&  \det\left(  \left(  \left(  x_{i}+y_{j}\right)  ^{m}\right)  _{1\leq i\leq
n,\ 1\leq j\leq n}\right) \nonumber\\
&  =\det\left(  \left(  P\left(  x_{i}+y_{j}\right)  \right)  _{1\leq i\leq
n,\ 1\leq j\leq n}\right) \nonumber\\
&  =\underbrace{p_{n-1}^{n}}_{\substack{=\delta_{n-1,m}^{n}\\\text{(since
}p_{n-1}=\delta_{n-1,m}\\\text{(by the definition of }p_{n-1}\text{))}%
}}\left(  \prod_{k=0}^{n-1}\dbinom{n-1}{k}\right)  \left(  \prod_{1\leq
i<j\leq n}\left(  x_{i}-x_{j}\right)  \right)  \left(  \prod_{1\leq i<j\leq
n}\left(  y_{j}-y_{i}\right)  \right) \nonumber\\
&  \ \ \ \ \ \ \ \ \ \ \left(  \text{by Exercise \ref{exe.vander-det.xi+yj}
\textbf{(c)}}\right) \nonumber\\
&  =\delta_{n-1,m}^{n}\left(  \prod_{k=0}^{n-1}\dbinom{n-1}{k}\right)  \left(
\prod_{1\leq i<j\leq n}\left(  x_{i}-x_{j}\right)  \right)  \left(
\prod_{1\leq i<j\leq n}\left(  y_{j}-y_{i}\right)  \right)  .
\label{sol.vander-det.xi+yj.a.almost}%
\end{align}

Now, let us forget that we fixed $m$. We thus have proven
(\ref{sol.vander-det.xi+yj.a.almost}) for every $m\in\left\{  0,1,\ldots
,n-1\right\}  $.

Now, let $m\in\left\{  0,1,\ldots,n-2\right\}  $. Thus, $m\neq n-1$, so that
$n-1\neq m$. Thus, $\delta_{n-1,m}=0$.

But $n$ is a positive integer. Hence, $0^{n}=0$. Now, taking the equality
$\delta_{n-1,m}=0$ to the $n$-th power, we obtain $\delta_{n-1,m}^{n}=0^{n}=0$.

On the other hand, $m\in\left\{  0,1,\ldots,n-2\right\}  \subseteq\left\{
0,1,\ldots,n-1\right\}  $, and therefore (\ref{sol.vander-det.xi+yj.a.almost})
holds. Hence,%
\begin{align*}
&  \det\left(  \left(  \left(  x_{i}+y_{j}\right)  ^{m}\right)  _{1\leq i\leq
n,\ 1\leq j\leq n}\right) \\
&  =\underbrace{\delta_{n-1,m}^{n}}_{=0}\left(  \prod_{k=0}^{n-1}\dbinom
{n-1}{k}\right)  \left(  \prod_{1\leq i<j\leq n}\left(  x_{i}-x_{j}\right)
\right)  \left(  \prod_{1\leq i<j\leq n}\left(  y_{j}-y_{i}\right)  \right) \\
&  =0.
\end{align*}
This solves Exercise \ref{exe.vander-det.xi+yj} \textbf{(a)}.

\textbf{(b)} For any two objects $i$ and $j$, we define $\delta_{i,j}$ to be
the element $%
\begin{cases}
1, & \text{if }i=j;\\
0, & \text{if }i\neq j
\end{cases}
$ of $\mathbb{K}$.

We have $\delta_{n-1,n-1}=1$ (since $n-1=n-1$) and thus $\delta_{n-1,n-1}%
^{n}=1^{n}=1$.

In our solution of Exercise \ref{exe.vander-det.xi+yj} \textbf{(a)}, we have
proven the equality (\ref{sol.vander-det.xi+yj.a.almost}) for every
$m\in\left\{  0,1,\ldots,n-1\right\}  $. Thus, we can apply this equality to
$m=n-1$. As a result, we obtain
\begin{align*}
&  \det\left(  \left(  \left(  x_{i}+y_{j}\right)  ^{n-1}\right)  _{1\leq
i\leq n,\ 1\leq j\leq n}\right) \\
&  =\underbrace{\delta_{n-1,n-1}^{n}}_{=1}\left(  \prod_{k=0}^{n-1}%
\dbinom{n-1}{k}\right)  \left(  \prod_{1\leq i<j\leq n}\left(  x_{i}%
-x_{j}\right)  \right)  \left(  \prod_{1\leq i<j\leq n}\left(  y_{j}%
-y_{i}\right)  \right) \\
&  =\left(  \prod_{k=0}^{n-1}\dbinom{n-1}{k}\right)  \left(  \prod_{1\leq
i<j\leq n}\left(  x_{i}-x_{j}\right)  \right)  \left(  \prod_{1\leq i<j\leq
n}\left(  y_{j}-y_{i}\right)  \right)  .
\end{align*}
This solves Exercise \ref{exe.vander-det.xi+yj} \textbf{(b)}.
\end{proof}
\end{verlong}

\subsection{Solution to Exercise \ref{exe.cauchy-det}}

We shall give a detailed solution to Exercise \ref{exe.cauchy-det} in Section
\ref{sect.sol.cauchy-det-lem} further below. For now, let us give some
pointers to the literature.

Exercise \ref{exe.cauchy-det} is obtained from \cite[Exercise 2.7.8(a)]%
{Reiner} by setting $a_{i}=x_{i}$ and $b_{j}=-y_{j}$. (See the ancillary PDF
file of the arXiv version of \cite{Reiner} for the solutions to the exercises.)

Exercise \ref{exe.cauchy-det} can also be obtained from \cite[Theorem
2]{Gri-19.9} by setting $k=\mathbb{K}$, $m=n$, $a_{j}=y_{j}$ and $b_{i}%
=-x_{i}$ (and observing that
\begin{align*}
&  \prod_{\substack{\left(  i,j\right)  \in\left\{  1,2,\ldots,n\right\}
^{2};\\i>j}}\underbrace{\left(  \left(  y_{i}-y_{j}\right)  \left(  \left(
-x_{j}\right)  -\left(  -x_{i}\right)  \right)  \right)  }_{\substack{=\left(
y_{j}-y_{i}\right)  \left(  x_{j}-x_{i}\right)  \\=\left(  x_{j}-x_{i}\right)
\left(  y_{j}-y_{i}\right)  }}\\
&  =\prod_{\substack{\left(  i,j\right)  \in\left\{  1,2,\ldots,n\right\}
^{2};\\i>j}}\left(  \left(  x_{j}-x_{i}\right)  \left(  y_{j}-y_{i}\right)
\right)  =\prod_{\substack{\left(  j,i\right)  \in\left\{  1,2,\ldots
,n\right\}  ^{2};\\j>i}}\left(  \left(  x_{i}-x_{j}\right)  \left(
y_{i}-y_{j}\right)  \right) \\
&  \ \ \ \ \ \ \ \ \ \ \left(  \text{here, we have renamed the index }\left(
i,j\right)  \text{ as }\left(  j,i\right)  \right) \\
&  =\prod_{1\leq i<j\leq n}\left(  \left(  x_{i}-x_{j}\right)  \left(
y_{i}-y_{j}\right)  \right)
\end{align*}
). The statement of \cite[Theorem 2]{Gri-19.9} makes the requirement that $k$
be a field; however, this is easily seen to be unnecessary for the proof.

\subsection{\label{sect.sol.cauchy-det-lem}Solution to Exercise
\ref{exe.cauchy-det-lem}}

\subsubsection{The solution}

Exercise \ref{exe.cauchy-det-lem} is the equality (13.70.3) in \cite[ancillary
PDF file]{Reiner}. The following solution is taken from \cite[ancillary PDF
file]{Reiner}.

\begin{proof}
[Solution to Exercise \ref{exe.cauchy-det-lem}.]We shall use the Iverson
bracket notation (introduced in Definition \ref{def.iverson}).

We have $n\in\left\{  1,2,\ldots,n\right\}  $ (since $n$ is a positive integer).

Define an $n\times n$-matrix $A\in\mathbb{K}^{n\times n}$ by $A=\left(
a_{i,j}\right)  _{1\leq i\leq n,\ 1\leq j\leq n}$.

For every $\left(  i,j\right)  \in\left\{  1,2,\ldots,n\right\}  ^{2}$, define
an element $b_{i,j}$ of $\mathbb{K}$ by
\begin{equation}
b_{i,j}=\left[  i=j\right]  a_{n,n}-\left[  i=n\text{ and }j\neq n\right]
a_{n,j}. \label{sol.cauchy-det-lem.beta.def}%
\end{equation}

Then, it is easy to see that
\begin{equation}
b_{i,j}=0\text{ for every }\left(  i,j\right)  \in\left\{  1,2,\ldots
,n\right\}  ^{2}\text{ satisfying }i<j. \label{sol.cauchy-det-lem.beta.triang}%
\end{equation}

[\textit{Proof of (\ref{sol.cauchy-det-lem.beta.triang}):} Let $\left(
i,j\right)  \in\left\{  1,2,\ldots,n\right\}  ^{2}$ be such that $i<j$. Then,
$i\neq j$ (since $i<j$), and thus we don't have $i=j$. Hence, $\left[
i=j\right]  =0$. Also, $j\in\left\{  1,2,\ldots,n\right\}  $ (since $\left(
i,j\right)  \in\left\{  1,2,\ldots,n\right\}  ^{2}$) and thus $j\leq n$, so
that $i<j\leq n$ and thus $i\neq n$. Hence, we don't have $i=n$. Thus, we
don't have $\left(  i=n\text{ and }j\neq n\right)  $ either. In other words,
$\left[  i=n\text{ and }j\neq n\right]  =0$. Now,
(\ref{sol.cauchy-det-lem.beta.def}) yields $b_{i,j}=\underbrace{\left[
i=j\right]  }_{=0}a_{n,n}-\underbrace{\left[  i=n\text{ and }j\neq n\right]
}_{=0}a_{n,j}=\underbrace{0a_{n,n}}_{=0}-\underbrace{0a_{n,j}}_{=0}=0$. This
proves (\ref{sol.cauchy-det-lem.beta.triang}).]

Define an $n\times n$-matrix $B\in\mathbb{K}^{n\times n}$ by $B=\left(
b_{i,j}\right)  _{1\leq i\leq n,\ 1\leq j\leq n}$. Thus, Exercise
\ref{exe.ps4.3} (applied to $B$ and $b_{i,j}$ instead of $A$ and $a_{i,j}$)
yields%
\begin{equation}
\det B=b_{1,1}b_{2,2}\cdots b_{n,n} \label{sol.cauchy-det-lem.detB=}%
\end{equation}
(since $b_{i,j}=0$ for every $\left(  i,j\right)  \in\left\{  1,2,\ldots
,n\right\}  ^{2}$ satisfying $i<j$).

But
\begin{equation}
\text{every }i\in\left\{  1,2,\ldots,n\right\}  \text{ satisfies }%
b_{i,i}=a_{n,n}. \label{sol.cauchy-det-lem.bii=}%
\end{equation}

[\textit{Proof of (\ref{sol.cauchy-det-lem.bii=}):} Let $i\in\left\{
1,2,\ldots,n\right\}  $. Then, we don't have $\left(  i=n\text{ and }i\neq
n\right)  $ (because the two statements $i=n$ and $i\neq n$ clearly contradict
each other). In other words, $\left[  i=n\text{ and }i\neq n\right]  =0$. Now,
(\ref{sol.cauchy-det-lem.beta.def}) (applied to $j=i$) yields%
\[
b_{i,i}=\underbrace{\left[  i=i\right]  }_{\substack{=1\\\text{(since
}i=i\text{ is true)}}}a_{n,n}-\underbrace{\left[  i=n\text{ and }i\neq
n\right]  }_{=0}a_{n,i}=\underbrace{1a_{n,n}}_{=a_{n,n}}-\underbrace{0a_{n,i}%
}_{=0}=a_{n,n}.
\]
This proves (\ref{sol.cauchy-det-lem.bii=}).]

Thus, (\ref{sol.cauchy-det-lem.detB=}) becomes
\begin{equation}
\det B=b_{1,1}b_{2,2}\cdots b_{n,n}=\prod_{i=1}^{n}\underbrace{b_{i,i}%
}_{\substack{=a_{n,n}\\\text{(by (\ref{sol.cauchy-det-lem.bii=}))}}%
}=\prod_{i=1}^{n}a_{n,n}=a_{n,n}^{n}. \label{sol.cauchy-det-lem.detB=2}%
\end{equation}

For each $\left(  i,j\right)  \in\left\{  1,2,\ldots,n\right\}  ^{2}$, we
define an element $c_{i,j}$ of $\mathbb{K}$ by%
\begin{equation}
c_{i,j}=\sum_{k=1}^{n}a_{i,k}b_{k,j}. \label{sol.cauchy-det-lem.cij=}%
\end{equation}

But we have $A=\left(  a_{i,j}\right)  _{1\leq i\leq n,\ 1\leq j\leq n}$ and
$B=\left(  b_{i,j}\right)  _{1\leq i\leq n,\ 1\leq j\leq n}$. The definition
of $AB$ thus yields
\[
AB=\left(  \underbrace{\sum_{k=1}^{n}a_{i,k}b_{k,j}}_{\substack{=c_{i,j}%
\\\text{(by (\ref{sol.cauchy-det-lem.cij=}))}}}\right)  _{1\leq i\leq
n,\ 1\leq j\leq n}=\left(  c_{i,j}\right)  _{1\leq i\leq n,\ 1\leq j\leq n}.
\]

But Theorem \ref{thm.det(AB)} yields
\begin{equation}
\det\left(  AB\right)  =\det A\cdot\underbrace{\det B}_{\substack{=a_{n,n}%
^{n}\\\text{(by (\ref{sol.cauchy-det-lem.detB=2}))}}}=\left(  \det A\right)
\cdot a_{n,n}^{n}=a_{n,n}^{n}\det A. \label{sol.cauchy-det-lem.detAB=}%
\end{equation}

But if $i$ and $j$ are any two elements of $\left\{  1,2,\ldots,n\right\}  $,
then%
\begin{equation}
c_{i,j}=a_{i,j}a_{n,n}-\left[  j\neq n\right]  a_{i,n}a_{n,j}.
\label{sol.cauchy-det-lem.cij=gen}%
\end{equation}

[\textit{Proof of (\ref{sol.cauchy-det-lem.cij=gen}):} Let $i$ and $j$ be two
elements of $\left\{  1,2,\ldots,n\right\}  $. The equality
(\ref{sol.cauchy-det-lem.cij=}) yields%
\begin{align}
c_{i,j}  &  =\underbrace{\sum_{k=1}^{n}}_{=\sum_{k\in\left\{  1,2,\ldots
,n\right\}  }}a_{i,k}\underbrace{b_{k,j}}_{\substack{=\left[  k=j\right]
a_{n,n}-\left[  k=n\text{ and }j\neq n\right]  a_{n,j}\\\text{(by
(\ref{sol.cauchy-det-lem.beta.def}) (applied to }k\text{ instead of
}i\text{))}}}\nonumber\\
&  =\sum_{k\in\left\{  1,2,\ldots,n\right\}  }\underbrace{a_{i,k}\left(
\left[  k=j\right]  a_{n,n}-\left[  k=n\text{ and }j\neq n\right]
a_{n,j}\right)  }_{=\left[  k=j\right]  a_{i,k}a_{n,n}-\left[  k=n\text{ and
}j\neq n\right]  a_{i,k}a_{n,j}}\nonumber\\
&  =\sum_{k\in\left\{  1,2,\ldots,n\right\}  }\left(  \left[  k=j\right]
a_{i,k}a_{n,n}-\left[  k=n\text{ and }j\neq n\right]  a_{i,k}a_{n,j}\right)
\nonumber\\
&  =\sum_{k\in\left\{  1,2,\ldots,n\right\}  }\left[  k=j\right]
a_{i,k}a_{n,n}-\sum_{k\in\left\{  1,2,\ldots,n\right\}  }\left[  k=n\text{ and
}j\neq n\right]  a_{i,k}a_{n,j}. \label{sol.cauchy-det-lem.cij=gen.pf.1}%
\end{align}
But
\begin{align*}
&  \sum_{k\in\left\{  1,2,\ldots,n\right\}  }\left[  k=j\right]
a_{i,k}a_{n,n}\\
&  =\underbrace{\left[  j=j\right]  }_{\substack{=1\\\text{(since }j=j\text{
holds)}}}a_{i,j}a_{n,n}+\sum_{\substack{k\in\left\{  1,2,\ldots,n\right\}
;\\k\neq j}}\underbrace{\left[  k=j\right]  }_{\substack{=0\\\text{(since
}k=j\text{ does not hold}\\\text{(since }k\neq j\text{))}}}a_{i,k}a_{n,n}\\
&  \ \ \ \ \ \ \ \ \ \ \left(
\begin{array}
[c]{c}%
\text{here, we have split off the addend for }k=j\text{ from the sum}\\
\text{(since }j\in\left\{  1,2,\ldots,n\right\}  \text{)}%
\end{array}
\right) \\
&  =a_{i,j}a_{n,n}+\underbrace{\sum_{\substack{k\in\left\{  1,2,\ldots
,n\right\}  ;\\k\neq j}}0a_{i,k}a_{n,n}}_{=0}=a_{i,j}a_{n,n}%
\end{align*}
and%
\begin{align*}
&  \sum_{k\in\left\{  1,2,\ldots,n\right\}  }\left[  k=n\text{ and }j\neq
n\right]  a_{i,k}a_{n,j}\\
&  =\underbrace{\left[  n=n\text{ and }j\neq n\right]  }_{\substack{=\left[
j\neq n\right]  \\\text{(since the statement }\left(  n=n\text{ and }j\neq
n\right)  \\\text{is equivalent to }\left(  j\neq n\right)  \\\text{(since
}n=n\text{ always holds))}}}a_{i,n}a_{n,j}+\sum_{\substack{k\in\left\{
1,2,\ldots,n\right\}  ;\\k\neq n}}\underbrace{\left[  k=n\text{ and }j\neq
n\right]  }_{\substack{=0\\\text{(since }\left(  k=n\text{ and }j\neq
n\right)  \text{ does not hold}\\\text{(since }k=n\text{ does not
hold}\\\text{(since }k\neq n\text{)))}}}a_{i,k}a_{n,j}\\
&  \ \ \ \ \ \ \ \ \ \ \left(
\begin{array}
[c]{c}%
\text{here, we have split off the addend for }k=n\text{ from the sum}\\
\text{(since }n\in\left\{  1,2,\ldots,n\right\}  \text{)}%
\end{array}
\right) \\
&  =\left[  j\neq n\right]  a_{i,n}a_{n,j}+\underbrace{\sum_{\substack{k\in
\left\{  1,2,\ldots,n\right\}  ;\\k\neq n}}0a_{i,k}a_{n,j}}_{=0}=\left[  j\neq
n\right]  a_{i,n}a_{n,j}.
\end{align*}

Hence, (\ref{sol.cauchy-det-lem.cij=gen.pf.1}) becomes
\begin{align*}
c_{i,j}  &  =\underbrace{\sum_{k\in\left\{  1,2,\ldots,n\right\}  }\left[
k=j\right]  a_{i,k}a_{n,n}}_{=a_{i,j}a_{n,n}}-\underbrace{\sum_{k\in\left\{
1,2,\ldots,n\right\}  }\left[  k=n\text{ and }j\neq n\right]  a_{i,k}a_{n,j}%
}_{=\left[  j\neq n\right]  a_{i,n}a_{n,j}}\\
&  =a_{i,j}a_{n,n}-\left[  j\neq n\right]  a_{i,n}a_{n,j}.
\end{align*}
This proves (\ref{sol.cauchy-det-lem.cij=gen}).]

Now, we can easily check that
\begin{equation}
c_{n,j}=0\text{ for every }j\in\left\{  1,2,\ldots,n-1\right\}  .
\label{sol.cauchy-det-lem.cnj=0}%
\end{equation}

[\textit{Proof of (\ref{sol.cauchy-det-lem.cnj=0}):} Let $j\in\left\{
1,2,\ldots,n-1\right\}  $. Then, $j\leq n-1<n$, so that $j\neq n$. Thus,
$\left[  j\neq n\right]  =1$. Also, $j\in\left\{  1,2,\ldots,n-1\right\}
\subseteq\left\{  1,2,\ldots,n\right\}  $ and $n\in\left\{  1,2,\ldots
,n\right\}  $. Hence, $n$ and $j$ are two elements of $\left\{  1,2,\ldots
,n\right\}  $. Thus, (\ref{sol.cauchy-det-lem.cij=gen}) (applied to $i=n$)
yields%
\[
c_{n,j}=\underbrace{a_{n,j}a_{n,n}}_{=a_{n,n}a_{n,j}}-\underbrace{\left[
j\neq n\right]  }_{=1}a_{n,n}a_{n,j}=a_{n,n}a_{n,j}-a_{n,n}a_{n,j}=0.
\]
This proves (\ref{sol.cauchy-det-lem.cnj=0}).]

So we have shown that $c_{n,j}=0$ for every $j\in\left\{  1,2,\ldots
,n-1\right\}  $. Hence, Theorem \ref{thm.laplace.pre} (applied to $AB$ and
$c_{i,j}$ instead of $A$ and $a_{i,j}$) yields%
\[
\det\left(  AB\right)  =c_{n,n}\cdot\det\left(  \left(  c_{i,j}\right)
_{1\leq i\leq n-1,\ 1\leq j\leq n-1}\right)
\]
(since $AB=\left(  c_{i,j}\right)  _{1\leq i\leq n,\ 1\leq j\leq n}$).
Comparing this with (\ref{sol.cauchy-det-lem.detAB=}), we find%
\begin{equation}
a_{n,n}^{n}\det A=c_{n,n}\cdot\det\left(  \left(  c_{i,j}\right)  _{1\leq
i\leq n-1,\ 1\leq j\leq n-1}\right)  . \label{sol.cauchy-det-lem.detA-and-c}%
\end{equation}

But every $\left(  i,j\right)  \in\left\{  1,2,\ldots,n-1\right\}  ^{2}$
satisfies%
\begin{equation}
c_{i,j}=a_{i,j}a_{n,n}-a_{i,n}a_{n,j}. \label{sol.cauchy-det-lem.cij=res}%
\end{equation}

[\textit{Proof of (\ref{sol.cauchy-det-lem.cij=res}):} Let $\left(
i,j\right)  \in\left\{  1,2,\ldots,n-1\right\}  ^{2}$. Since $\left(
i,j\right)  \in\left\{  1,2,\ldots,n-1\right\}  ^{2}$, we have $i\in\left\{
1,2,\ldots,n-1\right\}  $ and $j\in\left\{  1,2,\ldots,n-1\right\}  $. Since
$j\in\left\{  1,2,\ldots,n-1\right\}  $, we have $j\leq n-1<n$, so that $j\neq
n$. Thus, $\left[  j\neq n\right]  =1$.

From $i\in\left\{  1,2,\ldots,n-1\right\}  \subseteq\left\{  1,2,\ldots
,n\right\}  $ and $j\in\left\{  1,2,\ldots,n-1\right\}  \subseteq\left\{
1,2,\ldots,n\right\}  $, we conclude that $i$ and $j$ are two elements of
$\left\{  1,2,\ldots,n\right\}  $. Hence, (\ref{sol.cauchy-det-lem.cij=gen})
yields $c_{i,j}=a_{i,j}a_{n,n}-\underbrace{\left[  j\neq n\right]  }%
_{=1}a_{i,n}a_{n,j}=a_{i,j}a_{n,n}-a_{i,n}a_{n,j}$. This proves
(\ref{sol.cauchy-det-lem.cij=res}).]

Furthermore,
\begin{equation}
c_{n,n}=a_{n,n}^{2}. \label{sol.cauchy-det-lem.cnn=}%
\end{equation}

[\textit{Proof of (\ref{sol.cauchy-det-lem.cnn=}):} We have $n\in\left\{
1,2,\ldots,n\right\}  $. Thus, $n$ and $n$ are two elements of $\left\{
1,2,\ldots,n\right\}  $. Hence, (\ref{sol.cauchy-det-lem.cij=gen}) (applied to
$i=n$ and $j=n$) yields%
\[
c_{n,n}=a_{n,n}a_{n,n}-\underbrace{\left[  n\neq n\right]  }%
_{\substack{=0\\\text{(since we don't have }n\neq n\text{)}}}a_{n,n}%
a_{n,n}=a_{n,n}a_{n,n}-\underbrace{0a_{n,n}a_{n,n}}_{=0}=a_{n,n}%
a_{n,n}=a_{n,n}^{2}.
\]
This proves (\ref{sol.cauchy-det-lem.cnn=}).]

Now, (\ref{sol.cauchy-det-lem.detA-and-c}) becomes%
\begin{align*}
a_{n,n}^{n}\det A  &  =\underbrace{c_{n,n}}_{\substack{=a_{n,n}^{2}\\\text{(by
(\ref{sol.cauchy-det-lem.cnn=}))}}}\cdot\det\left(  \left(
\underbrace{c_{i,j}}_{\substack{=a_{i,j}a_{n,n}-a_{i,n}a_{n,j}\\\text{(by
(\ref{sol.cauchy-det-lem.cij=res}))}}}\right)  _{1\leq i\leq n-1,\ 1\leq j\leq
n-1}\right) \\
&  =a_{n,n}^{2}\cdot\det\left(  \left(  a_{i,j}a_{n,n}-a_{i,n}a_{n,j}\right)
_{1\leq i\leq n-1,\ 1\leq j\leq n-1}\right)  .
\end{align*}
We can divide both sides of this equality by $a_{n,n}^{2}$ (since $a_{n,n}%
^{2}$ is invertible in $\mathbb{K}$ (because $a_{n,n}$ is invertible in
$\mathbb{K}$)), and thus obtain%
\[
\dfrac{1}{a_{n,n}^{2}}a_{n,n}^{n}\det A=\det\left(  \left(  a_{i,j}%
a_{n,n}-a_{i,n}a_{n,j}\right)  _{1\leq i\leq n-1,\ 1\leq j\leq n-1}\right)  .
\]
Thus,%
\begin{align*}
\det\left(  \left(  a_{i,j}a_{n,n}-a_{i,n}a_{n,j}\right)  _{1\leq i\leq
n-1,\ 1\leq j\leq n-1}\right)   &  =\underbrace{\dfrac{1}{a_{n,n}^{2}}%
a_{n,n}^{n}}_{=a_{n,n}^{n-2}}\det\underbrace{A}_{=\left(  a_{i,j}\right)
_{1\leq i\leq n,\ 1\leq j\leq n}}\\
&  =a_{n,n}^{n-2}\cdot\det\left(  \left(  a_{i,j}\right)  _{1\leq i\leq
n,\ 1\leq j\leq n}\right)  .
\end{align*}
This solves Exercise \ref{exe.cauchy-det-lem}.
\end{proof}

\subsubsection{Solution to Exercise \ref{exe.cauchy-det}}

We shall now observe our promise and solve Exercise \ref{exe.cauchy-det} using
Exercise \ref{exe.cauchy-det-lem}. Again, our solution will follow
\cite[ancillary PDF file]{Reiner}.

We begin with some elementary lemmas:

\begin{lemma}
\label{lem.sol.cauchy-det-lem.1}Let $n\in\mathbb{N}$. Let $\left(
a_{i,j}\right)  _{1\leq i\leq n,\ 1\leq j\leq n}\in\mathbb{K}^{n\times n}$ be
an $n\times n$-matrix. Let $c_{1},c_{2},\ldots,c_{n}$ be $n$ elements of
$\mathbb{K}$. Then,
\[
\det\left(  \left(  c_{i}a_{i,j}\right)  _{1\leq i\leq n,\ 1\leq j\leq
n}\right)  =\left(  \prod_{i=1}^{n}c_{i}\right)  \cdot\det\left(  \left(
a_{i,j}\right)  _{1\leq i\leq n,\ 1\leq j\leq n}\right)  .
\]

\end{lemma}

\begin{proof}
[Proof of Lemma \ref{lem.sol.cauchy-det-lem.1}.]Define an $n\times n$-matrix
$A\in\mathbb{K}^{n\times n}$ by $A=\left(  a_{i,j}\right)  _{1\leq i\leq
n,\ 1\leq j\leq n}$.

Define an $n\times n$-matrix $C\in\mathbb{K}^{n\times n}$ by $C=\left(
c_{i}a_{i,j}\right)  _{1\leq i\leq n,\ 1\leq j\leq n}$. Hence,
(\ref{eq.det.eq.2}) (applied to $C$ and $c_{i}a_{i,j}$ instead of $A$ and
$a_{i,j}$) yields%
\begin{align*}
\det C  &  =\sum_{\sigma\in S_{n}}\left(  -1\right)  ^{\sigma}%
\underbrace{\prod_{i=1}^{n}\left(  c_{i}a_{i,\sigma\left(  i\right)  }\right)
}_{=\left(  \prod_{i=1}^{n}c_{i}\right)  \left(  \prod_{i=1}^{n}%
a_{i,\sigma\left(  i\right)  }\right)  }=\sum_{\sigma\in S_{n}}\left(
-1\right)  ^{\sigma}\left(  \prod_{i=1}^{n}c_{i}\right)  \left(  \prod
_{i=1}^{n}a_{i,\sigma\left(  i\right)  }\right) \\
&  =\left(  \prod_{i=1}^{n}c_{i}\right)  \cdot\underbrace{\sum_{\sigma\in
S_{n}}\left(  -1\right)  ^{\sigma}\prod_{i=1}^{n}a_{i,\sigma\left(  i\right)
}}_{\substack{=\det A\\\text{(by (\ref{eq.det.eq.2}))}}}=\left(  \prod
_{i=1}^{n}c_{i}\right)  \cdot\det\underbrace{A}_{=\left(  a_{i,j}\right)
_{1\leq i\leq n,\ 1\leq j\leq n}}\\
&  =\left(  \prod_{i=1}^{n}c_{i}\right)  \cdot\det\left(  \left(
a_{i,j}\right)  _{1\leq i\leq n,\ 1\leq j\leq n}\right)  .
\end{align*}
In view of $C=\left(  c_{i}a_{i,j}\right)  _{1\leq i\leq n,\ 1\leq j\leq n}$,
this rewrites as%
\[
\det\left(  \left(  c_{i}a_{i,j}\right)  _{1\leq i\leq n,\ 1\leq j\leq
n}\right)  =\left(  \prod_{i=1}^{n}c_{i}\right)  \cdot\det\left(  \left(
a_{i,j}\right)  _{1\leq i\leq n,\ 1\leq j\leq n}\right)  .
\]
This proves Lemma \ref{lem.sol.cauchy-det-lem.1}.
\end{proof}

\begin{lemma}
\label{lem.sol.cauchy-det-lem.2}Let $n\in\mathbb{N}$. Let $\left(
a_{i,j}\right)  _{1\leq i\leq n,\ 1\leq j\leq n}\in\mathbb{K}^{n\times n}$ be
an $n\times n$-matrix. Let $d_{1},d_{2},\ldots,d_{n}$ be $n$ elements of
$\mathbb{K}$. Then,
\[
\det\left(  \left(  d_{j}a_{i,j}\right)  _{1\leq i\leq n,\ 1\leq j\leq
n}\right)  =\left(  \prod_{j=1}^{n}d_{j}\right)  \cdot\det\left(  \left(
a_{i,j}\right)  _{1\leq i\leq n,\ 1\leq j\leq n}\right)  .
\]

\end{lemma}

\begin{proof}
[Proof of Lemma \ref{lem.sol.cauchy-det-lem.2}.]Let $\sigma\in S_{n}$. Thus,
$\sigma$ is a permutation of the set $\left\{  1,2,\ldots,n\right\}  $ (since
$S_{n}$ is the set of all permutations of the set $\left\{  1,2,\ldots
,n\right\}  $). In other words, $\sigma$ is a bijection $\left\{
1,2,\ldots,n\right\}  \rightarrow\left\{  1,2,\ldots,n\right\}  $. Now,%
\begin{align}
\underbrace{\prod_{i=1}^{n}}_{=\prod_{i\in\left\{  1,2,\ldots,n\right\}  }%
}d_{\sigma\left(  i\right)  }  &  =\prod_{i\in\left\{  1,2,\ldots,n\right\}
}d_{\sigma\left(  i\right)  }=\underbrace{\prod_{j\in\left\{  1,2,\ldots
,n\right\}  }}_{=\prod_{j=1}^{n}}d_{j}\nonumber\\
&  \ \ \ \ \ \ \ \ \ \ \left(
\begin{array}
[c]{c}%
\text{here, we have substituted }j\text{ for }\sigma\left(  i\right)  \text{
in the product,}\\
\text{since }\sigma\text{ is a bijection }\left\{  1,2,\ldots,n\right\}
\rightarrow\left\{  1,2,\ldots,n\right\}
\end{array}
\right) \nonumber\\
&  =\prod_{j=1}^{n}d_{j}. \label{pf.lem.sol.cauchy-det-lem.2.1}%
\end{align}

Now, forget that we fixed $\sigma$. We thus have proven the equality
(\ref{pf.lem.sol.cauchy-det-lem.2.1}) for each $\sigma\in S_{n}$.

Define an $n\times n$-matrix $A\in\mathbb{K}^{n\times n}$ by $A=\left(
a_{i,j}\right)  _{1\leq i\leq n,\ 1\leq j\leq n}$.

Define an $n\times n$-matrix $D\in\mathbb{K}^{n\times n}$ by $D=\left(
d_{j}a_{i,j}\right)  _{1\leq i\leq n,\ 1\leq j\leq n}$. Hence,
(\ref{eq.det.eq.2}) (applied to $D$ and $d_{j}a_{i,j}$ instead of $A$ and
$a_{i,j}$) yields%
\begin{align*}
\det D  &  =\sum_{\sigma\in S_{n}}\left(  -1\right)  ^{\sigma}%
\underbrace{\prod_{i=1}^{n}\left(  d_{\sigma\left(  i\right)  }a_{i,\sigma
\left(  i\right)  }\right)  }_{=\left(  \prod_{i=1}^{n}d_{\sigma\left(
i\right)  }\right)  \left(  \prod_{i=1}^{n}a_{i,\sigma\left(  i\right)
}\right)  }=\sum_{\sigma\in S_{n}}\left(  -1\right)  ^{\sigma}%
\underbrace{\left(  \prod_{i=1}^{n}d_{\sigma\left(  i\right)  }\right)
}_{\substack{=\prod_{j=1}^{n}d_{j}\\\text{(by
(\ref{pf.lem.sol.cauchy-det-lem.2.1}))}}}\left(  \prod_{i=1}^{n}%
a_{i,\sigma\left(  i\right)  }\right) \\
&  =\sum_{\sigma\in S_{n}}\left(  -1\right)  ^{\sigma}\left(  \prod_{j=1}%
^{n}d_{j}\right)  \left(  \prod_{i=1}^{n}a_{i,\sigma\left(  i\right)
}\right)  =\left(  \prod_{j=1}^{n}d_{j}\right)  \cdot\underbrace{\sum
_{\sigma\in S_{n}}\left(  -1\right)  ^{\sigma}\prod_{i=1}^{n}a_{i,\sigma
\left(  i\right)  }}_{\substack{=\det A\\\text{(by (\ref{eq.det.eq.2}))}}}\\
&  =\left(  \prod_{j=1}^{n}d_{j}\right)  \cdot\det\underbrace{A}_{=\left(
a_{i,j}\right)  _{1\leq i\leq n,\ 1\leq j\leq n}}=\left(  \prod_{j=1}^{n}%
d_{j}\right)  \cdot\det\left(  \left(  a_{i,j}\right)  _{1\leq i\leq n,\ 1\leq
j\leq n}\right)  .
\end{align*}
In view of $D=\left(  d_{j}a_{i,j}\right)  _{1\leq i\leq n,\ 1\leq j\leq n}$,
this rewrites as%
\[
\det\left(  \left(  d_{j}a_{i,j}\right)  _{1\leq i\leq n,\ 1\leq j\leq
n}\right)  =\left(  \prod_{j=1}^{n}d_{j}\right)  \cdot\det\left(  \left(
a_{i,j}\right)  _{1\leq i\leq n,\ 1\leq j\leq n}\right)  .
\]
This proves Lemma \ref{lem.sol.cauchy-det-lem.2}.
\end{proof}

Next, we will show a proposition that generalizes Exercise
\ref{exe.cauchy-det}:

\begin{proposition}
\label{prop.sol.cauchy-det.cauchy}Let $n\in\mathbb{N}$. For every
$i\in\left\{  1,2,\ldots,n\right\}  $, let $a_{i}$, $b_{i}$, $c_{i}$ and
$d_{i}$ be four elements of $\mathbb{K}$. Assume that $a_{i}d_{j}-b_{i}c_{j}$
is an invertible element of $\mathbb{K}$ for every $i\in\left\{
1,2,\ldots,n\right\}  $ and $j\in\left\{  1,2,\ldots,n\right\}  $. Then,
\[
\det\left(  \left(  \dfrac{1}{a_{i}d_{j}-b_{i}c_{j}}\right)  _{1\leq i\leq
n,\ 1\leq j\leq n}\right)  =\dfrac{\prod_{1\leq j<i\leq n}\left(  \left(
a_{i}b_{j}-a_{j}b_{i}\right)  \left(  c_{j}d_{i}-c_{i}d_{j}\right)  \right)
}{\prod_{\left(  i,j\right)  \in\left\{  1,2,\ldots,n\right\}  ^{2}}\left(
a_{i}d_{j}-b_{i}c_{j}\right)  }%
\]

\end{proposition}

\begin{proof}
[Proof of Proposition \ref{prop.sol.cauchy-det.cauchy}.]For every
$i\in\left\{  1,2,\ldots,n\right\}  $ and $j\in\left\{  1,2,\ldots,n\right\}
$, we define an element $x_{i,j}$ of $\mathbb{K}$ by%
\begin{equation}
x_{i,j}=a_{i}d_{j}-b_{i}c_{j}. \label{sol.cauchy-det.cauchy.xij=}%
\end{equation}

We have assumed that $a_{i}d_{j}-b_{i}c_{j}$ is an invertible element of
$\mathbb{K}$ for every $i\in\left\{  1,2,\ldots,n\right\}  $ and $j\in\left\{
1,2,\ldots,n\right\}  $. In other words, $x_{i,j}$ is an invertible element of
$\mathbb{K}$ for every $i\in\left\{  1,2,\ldots,n\right\}  $ and $j\in\left\{
1,2,\ldots,n\right\}  $ (because $x_{i,j}=a_{i}d_{j}-b_{i}c_{j}$ for every
$i\in\left\{  1,2,\ldots,n\right\}  $ and $j\in\left\{  1,2,\ldots,n\right\}
$).

We are going to show that
\begin{equation}
\det\left(  \left(  \dfrac{1}{x_{i,j}}\right)  _{1\leq i\leq k,\ 1\leq j\leq
k}\right)  =\dfrac{\prod_{1\leq j<i\leq k}\left(  \left(  a_{i}b_{j}%
-a_{j}b_{i}\right)  \left(  c_{j}d_{i}-c_{i}d_{j}\right)  \right)  }%
{\prod_{\left(  i,j\right)  \in\left\{  1,2,\ldots,k\right\}  ^{2}}x_{i,j}}
\label{sol.cauchy-det.cauchy.a.goal}%
\end{equation}
for every $k\in\left\{  0,1,\ldots,n\right\}  $.

[\textit{Proof of (\ref{sol.cauchy-det.cauchy.a.goal}):} We will prove
(\ref{sol.cauchy-det.cauchy.a.goal}) by induction over $k$:

\textit{Induction base:} It is easy to see that
(\ref{sol.cauchy-det.cauchy.a.goal}) holds if $k=0$%
\ \ \ \ \footnote{\textit{Proof.} Assume that $k=0$. Then, $\left(  \dfrac
{1}{x_{i,j}}\right)  _{1\leq i\leq k,\ 1\leq j\leq k}=\left(  \dfrac
{1}{x_{i,j}}\right)  _{1\leq i\leq0,\ 1\leq j\leq0}$ is a $0\times0$-matrix,
and thus has determinant $1$ (since every $0\times0$-matrix has determinant
$1$). In other words, $\det\left(  \left(  \dfrac{1}{x_{i,j}}\right)  _{1\leq
i\leq k,\ 1\leq j\leq k}\right)  =1$.
\par
But then there exist no integers $i$ and $j$ satisfying $1\leq j<i\leq k$
(because any such integers would have to satisfy $1\leq j<i\leq k=0<1$, which
is absurd). Thus, $\prod_{1\leq j<i\leq k}\left(  \left(  a_{i}b_{j}%
-a_{j}b_{i}\right)  \left(  c_{j}d_{i}-c_{i}d_{j}\right)  \right)  $ is an
empty product. Hence, $\prod_{1\leq j<i\leq k}\left(  \left(  a_{i}b_{j}%
-a_{j}b_{i}\right)  \left(  c_{j}d_{i}-c_{i}d_{j}\right)  \right)  =\left(
\text{empty product}\right)  =1$. Also, since $k=0$, we have $\left\{
1,2,\ldots,k\right\}  =\left\{  1,2,\ldots,0\right\}  =\varnothing$, so that
$\left\{  1,2,\ldots,k\right\}  ^{2}=\varnothing^{2}=\varnothing$. Thus,
$\prod_{\left(  i,j\right)  \in\left\{  1,2,\ldots,k\right\}  ^{2}}x_{i,j}$ is
an empty product. Hence, $\prod_{\left(  i,j\right)  \in\left\{
1,2,\ldots,k\right\}  ^{2}}x_{i,j}=\left(  \text{empty product}\right)  =1$.
\par
Now, dividing the equality $\prod_{1\leq j<i\leq k}\left(  \left(  a_{i}%
b_{j}-a_{j}b_{i}\right)  \left(  c_{j}d_{i}-c_{i}d_{j}\right)  \right)  =1$ by
the equality $\prod_{\left(  i,j\right)  \in\left\{  1,2,\ldots,k\right\}
^{2}}x_{i,j}=1$, we obtain%
\[
\dfrac{\prod_{1\leq j<i\leq k}\left(  \left(  a_{i}b_{j}-a_{j}b_{i}\right)
\left(  c_{j}d_{i}-c_{i}d_{j}\right)  \right)  }{\prod_{\left(  i,j\right)
\in\left\{  1,2,\ldots,k\right\}  ^{2}}x_{i,j}}=\dfrac{1}{1}=1.
\]
Comparing this with $\det\left(  \left(  \dfrac{1}{x_{i,j}}\right)  _{1\leq
i\leq k,\ 1\leq j\leq k}\right)  =1$, we obtain $\det\left(  \left(  \dfrac
{1}{x_{i,j}}\right)  _{1\leq i\leq k,\ 1\leq j\leq k}\right)  =\dfrac
{\prod_{1\leq j<i\leq k}\left(  \left(  a_{i}b_{j}-a_{j}b_{i}\right)  \left(
c_{j}d_{i}-c_{i}d_{j}\right)  \right)  }{\prod_{\left(  i,j\right)
\in\left\{  1,2,\ldots,k\right\}  ^{2}}x_{i,j}}$. Thus,
(\ref{sol.cauchy-det.cauchy.a.goal}) holds if $k=0$, qed.}. Hence, the
induction base is complete.

\textit{Induction step:} Let $m\in\left\{  1,2,\ldots,n\right\}  $. Assume
that (\ref{sol.cauchy-det.cauchy.a.goal}) is proven for $k=m-1$. We need to
show that (\ref{sol.cauchy-det.cauchy.a.goal}) holds for $k=m$.

We know that (\ref{sol.cauchy-det.cauchy.a.goal}) is proven for $k=m-1$. Thus,%
\begin{equation}
\det\left(  \left(  \dfrac{1}{x_{i,j}}\right)  _{1\leq i\leq m-1,\ 1\leq j\leq
m-1}\right)  =\dfrac{\prod_{1\leq j<i\leq m-1}\left(  \left(  a_{i}b_{j}%
-a_{j}b_{i}\right)  \left(  c_{j}d_{i}-c_{i}d_{j}\right)  \right)  }%
{\prod_{\left(  i,j\right)  \in\left\{  1,2,\ldots,m-1\right\}  ^{2}}x_{i,j}}.
\label{sol.cauchy-det.cauchy.a.indhyp}%
\end{equation}

We recall that $x_{i,j}$ is an invertible element of $\mathbb{K}$ for every
$i\in\left\{  1,2,\ldots,n\right\}  $ and $j\in\left\{  1,2,\ldots,n\right\}
$. Applying this to $i=m$ and $j=m$, we conclude that $x_{m,m}$ is an
invertible element of $\mathbb{K}$. Hence, $\dfrac{1}{x_{m,m}}$ is an
invertible element of $\mathbb{K}$ as well (because it is the inverse of
$x_{m,m}$). Thus, Exercise \ref{exe.cauchy-det-lem} (applied to $m$ and
$\dfrac{1}{x_{i,j}}$ instead of $n$ and $a_{i,j}$) yields
\begin{align*}
&  \det\left(  \left(  \dfrac{1}{x_{i,j}}\cdot\dfrac{1}{x_{m,m}}-\dfrac
{1}{x_{i,m}}\cdot\dfrac{1}{x_{m,j}}\right)  _{1\leq i\leq m-1,\ 1\leq j\leq
m-1}\right) \\
&  =\left(  \dfrac{1}{x_{m,m}}\right)  ^{m-2}\cdot\det\left(  \left(
\dfrac{1}{x_{i,j}}\right)  _{1\leq i\leq m,\ 1\leq j\leq m}\right)  .
\end{align*}
Multiplying both sides of this equality with $x_{m,m}^{m-2}$, we obtain
\begin{align}
&  x_{m,m}^{m-2}\cdot\det\left(  \left(  \dfrac{1}{x_{i,j}}\cdot\dfrac
{1}{x_{m,m}}-\dfrac{1}{x_{i,m}}\cdot\dfrac{1}{x_{m,j}}\right)  _{1\leq i\leq
m-1,\ 1\leq j\leq m-1}\right) \nonumber\\
&  =x_{m,m}^{m-2}\cdot\underbrace{\left(  \dfrac{1}{x_{m,m}}\right)  ^{m-2}%
}_{=\dfrac{1}{x_{m,m}^{m-2}}}\cdot\det\left(  \left(  \dfrac{1}{x_{i,j}%
}\right)  _{1\leq i\leq m,\ 1\leq j\leq m}\right) \nonumber\\
&  =\underbrace{x_{m,m}^{m-2}\cdot\dfrac{1}{x_{m,m}^{m-2}}}_{=1}\cdot
\det\left(  \left(  \dfrac{1}{x_{i,j}}\right)  _{1\leq i\leq m,\ 1\leq j\leq
m}\right) \nonumber\\
&  =\det\left(  \left(  \dfrac{1}{x_{i,j}}\right)  _{1\leq i\leq m,\ 1\leq
j\leq m}\right)  . \label{sol.cauchy-det.cauchy.a.1}%
\end{align}

But every $\left(  i,j\right)  \in\left\{  1,2,\ldots,m-1\right\}  ^{2}$
satisfies%
\begin{align}
&  \underbrace{x_{i,m}}_{\substack{=a_{i}d_{m}-b_{i}c_{m}\\\text{(by the
definition of }x_{i,m}\text{)}}}\underbrace{x_{m,j}}_{\substack{=a_{m}%
d_{j}-b_{m}c_{j}\\\text{(by the definition of }x_{m,j}\text{)}}%
}-\underbrace{x_{i,j}}_{\substack{=a_{i}d_{j}-b_{i}c_{j}\\\text{(by the
definition of }x_{i,j}\text{)}}}\underbrace{x_{m,m}}_{\substack{=a_{m}%
d_{m}-b_{m}c_{m}\\\text{(by the definition of }x_{m,m}\text{)}}}\nonumber\\
&  =\underbrace{\left(  a_{i}d_{m}-b_{i}c_{m}\right)  \left(  a_{m}d_{j}%
-b_{m}c_{j}\right)  }_{=a_{i}d_{m}a_{m}d_{j}-a_{i}d_{m}b_{m}c_{j}-b_{i}%
c_{m}a_{m}d_{j}+b_{i}c_{m}b_{m}c_{j}}-\underbrace{\left(  a_{i}d_{j}%
-b_{i}c_{j}\right)  \left(  a_{m}d_{m}-b_{m}c_{m}\right)  }_{=a_{i}d_{j}%
a_{m}d_{m}-a_{i}d_{j}b_{m}c_{m}-b_{i}c_{j}a_{m}d_{m}+b_{i}c_{j}b_{m}c_{m}%
}\nonumber\\
&  =\left(  \underbrace{a_{i}d_{m}a_{m}d_{j}}_{=a_{i}d_{j}a_{m}d_{m}%
}-\underbrace{a_{i}d_{m}b_{m}c_{j}}_{=a_{i}b_{m}c_{j}d_{m}}-\underbrace{b_{i}%
c_{m}a_{m}d_{j}}_{=a_{m}b_{i}c_{m}d_{j}}+\underbrace{b_{i}c_{m}b_{m}c_{j}%
}_{=b_{i}c_{j}b_{m}c_{m}}\right) \nonumber\\
&  \ \ \ \ \ \ \ \ \ \ -\left(  a_{i}d_{j}a_{m}d_{m}-\underbrace{a_{i}%
d_{j}b_{m}c_{m}}_{=a_{i}b_{m}c_{m}d_{j}}-\underbrace{b_{i}c_{j}a_{m}d_{m}%
}_{=a_{m}b_{i}c_{j}d_{m}}+b_{i}c_{j}b_{m}c_{m}\right) \nonumber\\
&  =\left(  a_{i}d_{j}a_{m}d_{m}-a_{i}b_{m}c_{j}d_{m}-a_{m}b_{i}c_{m}%
d_{j}+b_{i}c_{j}b_{m}c_{m}\right) \nonumber\\
&  \ \ \ \ \ \ \ \ \ \ -\left(  a_{i}d_{j}a_{m}d_{m}-a_{i}b_{m}c_{m}%
d_{j}-a_{m}b_{i}c_{j}d_{m}+b_{i}c_{j}b_{m}c_{m}\right) \nonumber\\
&  =a_{m}b_{i}c_{j}d_{m}-a_{m}b_{i}c_{m}d_{j}-a_{i}b_{m}c_{j}d_{m}+a_{i}%
b_{m}c_{m}d_{j}\nonumber\\
&  =\left(  a_{m}b_{i}-a_{i}b_{m}\right)  \left(  c_{j}d_{m}-c_{m}%
d_{j}\right)  . \label{sol.cauchy-det.cauchy.a.diff-fact}%
\end{align}
Hence, every $\left(  i,j\right)  \in\left\{  1,2,\ldots,m-1\right\}  ^{2}$
satisfies
\begin{align*}
&  \underbrace{\dfrac{1}{x_{i,j}}\cdot\dfrac{1}{x_{m,m}}}_{=\dfrac{1}%
{x_{i,j}x_{m,m}}}-\underbrace{\dfrac{1}{x_{i,m}}\cdot\dfrac{1}{x_{m,j}}%
}_{=\dfrac{1}{x_{i,m}x_{m,j}}}\\
&  =\dfrac{1}{x_{i,j}x_{m,m}}-\dfrac{1}{x_{i,m}x_{m,j}}=\dfrac{x_{i,m}%
x_{m,j}-x_{i,j}x_{m,m}}{x_{i,j}x_{m,m}x_{i,m}x_{m,j}}\\
&  =\dfrac{1}{x_{i,j}x_{m,m}x_{i,m}x_{m,j}}\cdot\underbrace{\left(
x_{i,m}x_{m,j}-x_{i,j}x_{m,m}\right)  }_{\substack{=\left(  a_{m}b_{i}%
-a_{i}b_{m}\right)  \left(  c_{j}d_{m}-c_{m}d_{j}\right)  \\\text{(by
(\ref{sol.cauchy-det.cauchy.a.diff-fact}))}}}\\
&  =\dfrac{1}{x_{i,j}x_{m,m}x_{i,m}x_{m,j}}\cdot\left(  a_{m}b_{i}-a_{i}%
b_{m}\right)  \left(  c_{j}d_{m}-c_{m}d_{j}\right) \\
&  =\dfrac{c_{j}d_{m}-c_{m}d_{j}}{x_{m,m}x_{m,j}}\cdot\dfrac{a_{m}b_{i}%
-a_{i}b_{m}}{x_{i,m}}\cdot\dfrac{1}{x_{i,j}}.
\end{align*}
Hence,
\begin{align*}
&  \det\left(  \left(  \underbrace{\dfrac{1}{x_{i,j}}\cdot\dfrac{1}{x_{m,m}%
}-\dfrac{1}{x_{i,m}}\cdot\dfrac{1}{x_{m,j}}}_{=\dfrac{c_{j}d_{m}-c_{m}d_{j}%
}{x_{m,m}x_{m,j}}\cdot\dfrac{a_{m}b_{i}-a_{i}b_{m}}{x_{i,m}}\cdot\dfrac
{1}{x_{i,j}}}\right)  _{1\leq i\leq m-1,\ 1\leq j\leq m-1}\right) \\
&  =\det\left(  \left(  \dfrac{c_{j}d_{m}-c_{m}d_{j}}{x_{m,m}x_{m,j}}%
\cdot\dfrac{a_{m}b_{i}-a_{i}b_{m}}{x_{i,m}}\cdot\dfrac{1}{x_{i,j}}\right)
_{1\leq i\leq m-1,\ 1\leq j\leq m-1}\right) \\
&  =\left(  \prod_{j=1}^{m-1}\underbrace{\dfrac{c_{j}d_{m}-c_{m}d_{j}}%
{x_{m,m}x_{m,j}}}_{=\dfrac{1}{x_{m,m}}\cdot\dfrac{c_{j}d_{m}-c_{m}d_{j}%
}{x_{m,j}}}\right)  \cdot\underbrace{\det\left(  \left(  \dfrac{a_{m}%
b_{i}-a_{i}b_{m}}{x_{i,m}}\cdot\dfrac{1}{x_{i,j}}\right)  _{1\leq i\leq
m-1,\ 1\leq j\leq m-1}\right)  }_{\substack{=\left(  \prod_{i=1}^{m-1}%
\dfrac{a_{m}b_{i}-a_{i}b_{m}}{x_{i,m}}\right)  \cdot\det\left(  \left(
\dfrac{1}{x_{i,j}}\right)  _{1\leq i\leq m-1,\ 1\leq j\leq m-1}\right)
\\\text{(by Lemma \ref{lem.sol.cauchy-det-lem.1} (applied to }m-1\text{,
}\dfrac{1}{x_{i,j}}\\\text{and }\dfrac{a_{m}b_{i}-a_{i}b_{m}}{x_{i,m}}\text{
instead of }n\text{, }a_{i,j}\text{ and }c_{i}\text{))}}}\\
&  \ \ \ \ \ \ \ \ \ \ \left(
\begin{array}
[c]{c}%
\text{by Lemma \ref{lem.sol.cauchy-det-lem.2} (applied to }m-1\text{, }%
\dfrac{a_{m}b_{i}-a_{i}b_{m}}{x_{i,m}}\cdot\dfrac{1}{x_{i,j}}\text{ and}\\
\dfrac{c_{j}d_{m}-c_{m}d_{j}}{x_{m,m}x_{m,j}}\text{ instead of }n\text{,
}a_{i,j}\text{ and }d_{j}\text{)}%
\end{array}
\right) \\
&  =\underbrace{\left(  \prod_{j=1}^{m-1}\left(  \dfrac{1}{x_{m,m}}\cdot
\dfrac{c_{j}d_{m}-c_{m}d_{j}}{x_{m,j}}\right)  \right)  }_{=\left(
\prod_{j=1}^{m-1}\dfrac{1}{x_{m,m}}\right)  \cdot\left(  \prod_{j=1}%
^{m-1}\dfrac{c_{j}d_{m}-c_{m}d_{j}}{x_{m,j}}\right)  }\cdot\left(  \prod
_{i=1}^{m-1}\dfrac{a_{m}b_{i}-a_{i}b_{m}}{x_{i,m}}\right) \\
&  \ \ \ \ \ \ \ \ \ \ \cdot\underbrace{\det\left(  \left(  \dfrac{1}{x_{i,j}%
}\right)  _{1\leq i\leq m-1,\ 1\leq j\leq m-1}\right)  }_{\substack{=\dfrac
{\prod_{1\leq j<i\leq m-1}\left(  \left(  a_{i}b_{j}-a_{j}b_{i}\right)
\left(  c_{j}d_{i}-c_{i}d_{j}\right)  \right)  }{\prod_{\left(  i,j\right)
\in\left\{  1,2,\ldots,m-1\right\}  ^{2}}x_{i,j}}\\\text{(by
(\ref{sol.cauchy-det.cauchy.a.indhyp}))}}}
\end{align*}%
\begin{align}
&  =\underbrace{\left(  \prod_{j=1}^{m-1}\dfrac{1}{x_{m,m}}\right)
}_{=\left(  \dfrac{1}{x_{m,m}}\right)  ^{m-1}}\cdot\left(  \prod_{j=1}%
^{m-1}\underbrace{\dfrac{c_{j}d_{m}-c_{m}d_{j}}{x_{m,j}}}_{=\dfrac{1}{x_{m,j}%
}\cdot\left(  c_{j}d_{m}-c_{m}d_{j}\right)  }\right)  \cdot\left(  \prod
_{i=1}^{m-1}\underbrace{\dfrac{a_{m}b_{i}-a_{i}b_{m}}{x_{i,m}}}_{=\dfrac
{1}{x_{i,m}}\cdot\left(  a_{m}b_{i}-a_{i}b_{m}\right)  }\right) \nonumber\\
&  \ \ \ \ \ \ \ \ \ \ \cdot\underbrace{\dfrac{\prod_{1\leq j<i\leq
m-1}\left(  \left(  a_{i}b_{j}-a_{j}b_{i}\right)  \left(  c_{j}d_{i}%
-c_{i}d_{j}\right)  \right)  }{\prod_{\left(  i,j\right)  \in\left\{
1,2,\ldots,m-1\right\}  ^{2}}x_{i,j}}}_{=\dfrac{1}{\prod_{\left(  i,j\right)
\in\left\{  1,2,\ldots,m-1\right\}  ^{2}}x_{i,j}}\cdot\left(  \prod_{1\leq
j<i\leq m-1}\left(  \left(  a_{i}b_{j}-a_{j}b_{i}\right)  \left(  c_{j}%
d_{i}-c_{i}d_{j}\right)  \right)  \right)  }\nonumber\\
&  =\left(  \dfrac{1}{x_{m,m}}\right)  ^{m-1}\cdot\underbrace{\left(
\prod_{j=1}^{m-1}\left(  \dfrac{1}{x_{m,j}}\cdot\left(  c_{j}d_{m}-c_{m}%
d_{j}\right)  \right)  \right)  }_{=\left(  \prod_{j=1}^{m-1}\dfrac{1}%
{x_{m,j}}\right)  \cdot\left(  \prod_{j=1}^{m-1}\left(  c_{j}d_{m}-c_{m}%
d_{j}\right)  \right)  }\cdot\underbrace{\left(  \prod_{i=1}^{m-1}\left(
\dfrac{1}{x_{i,m}}\cdot\left(  a_{m}b_{i}-a_{i}b_{m}\right)  \right)  \right)
}_{=\left(  \prod_{i=1}^{m-1}\dfrac{1}{x_{i,m}}\right)  \cdot\left(
\prod_{i=1}^{m-1}\left(  a_{m}b_{i}-a_{i}b_{m}\right)  \right)  }\nonumber\\
&  \ \ \ \ \ \ \ \ \ \ \cdot\underbrace{\dfrac{1}{\prod_{\left(  i,j\right)
\in\left\{  1,2,\ldots,m-1\right\}  ^{2}}x_{i,j}}}_{=\prod_{\left(
i,j\right)  \in\left\{  1,2,\ldots,m-1\right\}  ^{2}}\dfrac{1}{x_{i,j}}}%
\cdot\underbrace{\left(  \prod_{1\leq j<i\leq m-1}\left(  \left(  a_{i}%
b_{j}-a_{j}b_{i}\right)  \left(  c_{j}d_{i}-c_{i}d_{j}\right)  \right)
\right)  }_{=\left(  \prod_{1\leq j<i\leq m-1}\left(  a_{i}b_{j}-a_{j}%
b_{i}\right)  \right)  \left(  \prod_{1\leq j<i\leq m-1}\left(  c_{j}%
d_{i}-c_{i}d_{j}\right)  \right)  }\nonumber\\
&  =\left(  \dfrac{1}{x_{m,m}}\right)  ^{m-1}\cdot\left(  \prod_{j=1}%
^{m-1}\dfrac{1}{x_{m,j}}\right)  \cdot\left(  \prod_{j=1}^{m-1}\left(
c_{j}d_{m}-c_{m}d_{j}\right)  \right)  \cdot\left(  \prod_{i=1}^{m-1}\dfrac
{1}{x_{i,m}}\right)  \cdot\left(  \prod_{i=1}^{m-1}\left(  a_{m}b_{i}%
-a_{i}b_{m}\right)  \right) \nonumber\\
&  \ \ \ \ \ \ \ \ \ \ \cdot\left(  \prod_{\left(  i,j\right)  \in\left\{
1,2,\ldots,m-1\right\}  ^{2}}\dfrac{1}{x_{i,j}}\right)  \cdot\left(
\prod_{1\leq j<i\leq m-1}\left(  a_{i}b_{j}-a_{j}b_{i}\right)  \right)
\left(  \prod_{1\leq j<i\leq m-1}\left(  c_{j}d_{i}-c_{i}d_{j}\right)  \right)
\nonumber\\
&  =\left(  \dfrac{1}{x_{m,m}}\right)  ^{m-1}\cdot\left(  \prod_{j=1}%
^{m-1}\dfrac{1}{x_{m,j}}\right)  \cdot\left(  \prod_{i=1}^{m-1}\dfrac
{1}{x_{i,m}}\right)  \cdot\left(  \prod_{\left(  i,j\right)  \in\left\{
1,2,\ldots,m-1\right\}  ^{2}}\dfrac{1}{x_{i,j}}\right) \nonumber\\
&  \ \ \ \ \ \ \ \ \ \ \cdot\left(  \prod_{i=1}^{m-1}\left(  a_{m}b_{i}%
-a_{i}b_{m}\right)  \right)  \cdot\left(  \prod_{1\leq j<i\leq m-1}\left(
a_{i}b_{j}-a_{j}b_{i}\right)  \right) \nonumber\\
&  \ \ \ \ \ \ \ \ \ \ \cdot\left(  \prod_{j=1}^{m-1}\left(  c_{j}d_{m}%
-c_{m}d_{j}\right)  \right)  \cdot\left(  \prod_{1\leq j<i\leq m-1}\left(
c_{j}d_{i}-c_{i}d_{j}\right)  \right)  . \label{sol.cauchy-det.cauchy.a.3}%
\end{align}
But%
\begin{align*}
&  \underbrace{\prod_{\left(  i,j\right)  \in\left\{  1,2,\ldots,m\right\}
^{2}}}_{=\prod_{i\in\left\{  1,2,\ldots,m\right\}  }\prod_{j\in\left\{
1,2,\ldots,m\right\}  }}x_{i,j}\\
&  =\underbrace{\prod_{i\in\left\{  1,2,\ldots,m\right\}  }}_{=\prod_{i=1}%
^{m}}\underbrace{\prod_{j\in\left\{  1,2,\ldots,m\right\}  }}_{=\prod
_{j=1}^{m}}x_{i,j}=\prod_{i=1}^{m}\underbrace{\prod_{j=1}^{m}x_{i,j}%
}_{\substack{=\left(  \prod_{j=1}^{m-1}x_{i,j}\right)  \cdot x_{i,m}%
\\\text{(here, we have split off the}\\\text{factor for }j=m\text{ from the
product)}}}\\
&  =\prod_{i=1}^{m}\left(  \left(  \prod_{j=1}^{m-1}x_{i,j}\right)  \cdot
x_{i,m}\right) \\
&  =\underbrace{\left(  \prod_{i=1}^{m}\prod_{j=1}^{m-1}x_{i,j}\right)
}_{\substack{=\left(  \prod_{i=1}^{m-1}\prod_{j=1}^{m-1}x_{i,j}\right)
\cdot\left(  \prod_{j=1}^{m-1}x_{m,j}\right)  \\\text{(here, we have split off
the factor for }i=m\text{ from the outer product)}}}\cdot\underbrace{\left(
\prod_{i=1}^{m}x_{i,m}\right)  }_{\substack{=\left(  \prod_{i=1}^{m-1}%
x_{i,m}\right)  \cdot x_{m,m}\\\text{(here, we have split off the factor for
}i=m\\\text{from the product)}}}\\
&  =\left(  \underbrace{\prod_{i=1}^{m-1}}_{=\prod_{i\in\left\{
1,2,\ldots,m-1\right\}  }}\underbrace{\prod_{j=1}^{m-1}}_{=\prod_{j\in\left\{
1,2,\ldots,m-1\right\}  }}x_{i,j}\right)  \cdot\left(  \prod_{j=1}%
^{m-1}x_{m,j}\right)  \cdot\left(  \prod_{i=1}^{m-1}x_{i,m}\right)  \cdot
x_{m,m}\\
&  =\left(  \underbrace{\prod_{i\in\left\{  1,2,\ldots,m-1\right\}  }%
\prod_{j\in\left\{  1,2,\ldots,m-1\right\}  }}_{=\prod_{\left(  i,j\right)
\in\left\{  1,2,\ldots,m-1\right\}  ^{2}}}x_{i,j}\right)  \cdot\left(
\prod_{j=1}^{m-1}x_{m,j}\right)  \cdot\left(  \prod_{i=1}^{m-1}x_{i,m}\right)
\cdot x_{m,m}\\
&  =\left(  \prod_{\left(  i,j\right)  \in\left\{  1,2,\ldots,m-1\right\}
^{2}}x_{i,j}\right)  \cdot\left(  \prod_{j=1}^{m-1}x_{m,j}\right)
\cdot\left(  \prod_{i=1}^{m-1}x_{i,m}\right)  \cdot x_{m,m}.
\end{align*}
Taking the (multiplicative) inverse of both sides of this equality, we obtain%
\begin{align*}
&  \dfrac{1}{\prod_{\left(  i,j\right)  \in\left\{  1,2,\ldots,m\right\}
^{2}}x_{i,j}}\\
&  =\dfrac{1}{\left(  \prod_{\left(  i,j\right)  \in\left\{  1,2,\ldots
,m-1\right\}  ^{2}}x_{i,j}\right)  \cdot\left(  \prod_{j=1}^{m-1}%
x_{m,j}\right)  \cdot\left(  \prod_{i=1}^{m-1}x_{i,m}\right)  \cdot x_{m,m}}\\
&  =\underbrace{\dfrac{1}{\prod_{\left(  i,j\right)  \in\left\{
1,2,\ldots,m-1\right\}  ^{2}}x_{i,j}}}_{=\prod_{\left(  i,j\right)
\in\left\{  1,2,\ldots,m-1\right\}  ^{2}}\dfrac{1}{x_{i,j}}}\cdot
\underbrace{\dfrac{1}{\prod_{j=1}^{m-1}x_{m,j}}}_{=\prod_{j=1}^{m-1}\dfrac
{1}{x_{m,j}}}\cdot\underbrace{\dfrac{1}{\prod_{i=1}^{m-1}x_{i,m}}}%
_{=\prod_{i=1}^{m-1}\dfrac{1}{x_{i,m}}}\cdot\dfrac{1}{x_{m,m}}\\
&  =\left(  \prod_{\left(  i,j\right)  \in\left\{  1,2,\ldots,m-1\right\}
^{2}}\dfrac{1}{x_{i,j}}\right)  \cdot\left(  \prod_{j=1}^{m-1}\dfrac
{1}{x_{m,j}}\right)  \cdot\left(  \prod_{i=1}^{m-1}\dfrac{1}{x_{i,m}}\right)
\cdot\dfrac{1}{x_{m,m}}\\
&  =\dfrac{1}{x_{m,m}}\cdot\left(  \prod_{j=1}^{m-1}\dfrac{1}{x_{m,j}}\right)
\cdot\left(  \prod_{i=1}^{m-1}\dfrac{1}{x_{i,m}}\right)  \cdot\left(
\prod_{\left(  i,j\right)  \in\left\{  1,2,\ldots,m-1\right\}  ^{2}}\dfrac
{1}{x_{i,j}}\right)  .
\end{align*}
Thus,
\begin{align}
&  \dfrac{1}{x_{m,m}}\cdot\left(  \prod_{j=1}^{m-1}\dfrac{1}{x_{m,j}}\right)
\cdot\left(  \prod_{i=1}^{m-1}\dfrac{1}{x_{i,m}}\right)  \cdot\left(
\prod_{\left(  i,j\right)  \in\left\{  1,2,\ldots,m-1\right\}  ^{2}}\dfrac
{1}{x_{i,j}}\right) \nonumber\\
&  =\dfrac{1}{\prod_{\left(  i,j\right)  \in\left\{  1,2,\ldots,m\right\}
^{2}}x_{i,j}}. \label{sol.cauchy-det.cauchy.a.5a}%
\end{align}
Multiplying both sides of this equality by $x_{m,m}$, we obtain%
\begin{align}
&  \left(  \prod_{j=1}^{m-1}\dfrac{1}{x_{m,j}}\right)  \cdot\left(
\prod_{i=1}^{m-1}\dfrac{1}{x_{i,m}}\right)  \cdot\left(  \prod_{\left(
i,j\right)  \in\left\{  1,2,\ldots,m-1\right\}  ^{2}}\dfrac{1}{x_{i,j}}\right)
\nonumber\\
&  =x_{m,m}\cdot\dfrac{1}{\prod_{\left(  i,j\right)  \in\left\{
1,2,\ldots,m\right\}  ^{2}}x_{i,j}}. \label{sol.cauchy-det.cauchy.a.5aa}%
\end{align}

On the other hand,
\begin{align}
&  \underbrace{\prod_{1\leq j<i\leq m}}_{=\prod_{i=1}^{m}\prod_{j=1}^{i-1}%
}\left(  a_{i}b_{j}-a_{j}b_{i}\right) \nonumber\\
&  =\prod_{i=1}^{m}\prod_{j=1}^{i-1}\left(  a_{i}b_{j}-a_{j}b_{i}\right)
=\left(  \underbrace{\prod_{i=1}^{m-1}\prod_{j=1}^{i-1}}_{=\prod_{1\leq
j<i\leq m-1}}\left(  a_{i}b_{j}-a_{j}b_{i}\right)  \right)  \cdot\prod
_{j=1}^{m-1}\left(  a_{m}b_{j}-a_{j}b_{m}\right) \nonumber\\
&  \ \ \ \ \ \ \ \ \ \ \left(  \text{here, we have split off the factor for
}i=m\text{ from the outer product}\right) \nonumber\\
&  =\left(  \prod_{1\leq j<i\leq m-1}\left(  a_{i}b_{j}-a_{j}b_{i}\right)
\right)  \cdot\underbrace{\left(  \prod_{j=1}^{m-1}\left(  a_{m}b_{j}%
-a_{j}b_{m}\right)  \right)  }_{\substack{=\prod_{i=1}^{m-1}\left(  a_{m}%
b_{i}-a_{i}b_{m}\right)  \\\text{(here, we renamed the index }j\\\text{as
}i\text{ in the product)}}}\nonumber\\
&  =\left(  \prod_{1\leq j<i\leq m-1}\left(  a_{i}b_{j}-a_{j}b_{i}\right)
\right)  \cdot\left(  \prod_{i=1}^{m-1}\left(  a_{m}b_{i}-a_{i}b_{m}\right)
\right) \nonumber\\
&  =\left(  \prod_{i=1}^{m-1}\left(  a_{m}b_{i}-a_{i}b_{m}\right)  \right)
\cdot\left(  \prod_{1\leq j<i\leq m-1}\left(  a_{i}b_{j}-a_{j}b_{i}\right)
\right)  . \label{sol.cauchy-det.cauchy.a.5b}%
\end{align}

Also,%
\begin{align}
&  \underbrace{\prod_{1\leq j<i\leq m}}_{=\prod_{i=1}^{m}\prod_{j=1}^{i-1}%
}\left(  c_{j}d_{i}-c_{i}d_{j}\right) \nonumber\\
&  =\prod_{i=1}^{m}\prod_{j=1}^{i-1}\left(  c_{j}d_{i}-c_{i}d_{j}\right)
=\left(  \underbrace{\prod_{i=1}^{m-1}\prod_{j=1}^{i-1}}_{=\prod_{1\leq
j<i\leq m-1}}\left(  c_{j}d_{i}-c_{i}d_{j}\right)  \right)  \cdot\prod
_{j=1}^{m-1}\left(  c_{j}d_{m}-c_{m}d_{j}\right) \nonumber\\
&  \ \ \ \ \ \ \ \ \ \ \left(  \text{here, we have split off the factor for
}i=m\text{ from the outer product}\right) \nonumber\\
&  =\left(  \prod_{1\leq j<i\leq m-1}\left(  c_{j}d_{i}-c_{i}d_{j}\right)
\right)  \cdot\prod_{j=1}^{m-1}\left(  c_{j}d_{m}-c_{m}d_{j}\right)
\nonumber\\
&  =\left(  \prod_{j=1}^{m-1}\left(  c_{j}d_{m}-c_{m}d_{j}\right)  \right)
\cdot\left(  \prod_{1\leq j<i\leq m-1}\left(  c_{j}d_{i}-c_{i}d_{j}\right)
\right)  . \label{sol.cauchy-det.cauchy.a.5c}%
\end{align}
Now, (\ref{sol.cauchy-det.cauchy.a.3}) becomes%
\begin{align*}
&  \det\left(  \left(  \dfrac{1}{x_{i,j}}\cdot\dfrac{1}{x_{m,m}}-\dfrac
{1}{x_{i,m}}\cdot\dfrac{1}{x_{m,j}}\right)  _{1\leq i\leq m-1,\ 1\leq j\leq
m-1}\right) \\
&  =\underbrace{\left(  \dfrac{1}{x_{m,m}}\right)  ^{m-1}}_{\substack{=\left(
\dfrac{1}{x_{m,m}}\right)  ^{\left(  m-2\right)  +1}\\=\left(  \dfrac
{1}{x_{m,m}}\right)  ^{m-2}\left(  \dfrac{1}{x_{m,m}}\right)  }}\cdot
\underbrace{\left(  \prod_{j=1}^{m-1}\dfrac{1}{x_{m,j}}\right)  \cdot\left(
\prod_{i=1}^{m-1}\dfrac{1}{x_{i,m}}\right)  \cdot\left(  \prod_{\left(
i,j\right)  \in\left\{  1,2,\ldots,m-1\right\}  ^{2}}\dfrac{1}{x_{i,j}%
}\right)  }_{\substack{=x_{m,m}\cdot\dfrac{1}{\prod_{\left(  i,j\right)
\in\left\{  1,2,\ldots,m\right\}  ^{2}}x_{i,j}}\\\text{(by
(\ref{sol.cauchy-det.cauchy.a.5aa}))}}}\\
&  \ \ \ \ \ \ \ \ \ \ \cdot\underbrace{\left(  \prod_{i=1}^{m-1}\left(
a_{m}b_{i}-a_{i}b_{m}\right)  \right)  \cdot\left(  \prod_{1\leq j<i\leq
m-1}\left(  a_{i}b_{j}-a_{j}b_{i}\right)  \right)  }_{\substack{=\prod_{1\leq
j<i\leq m}\left(  a_{i}b_{j}-a_{j}b_{i}\right)  \\\text{(by
(\ref{sol.cauchy-det.cauchy.a.5b}))}}}\\
&  \ \ \ \ \ \ \ \ \ \ \cdot\underbrace{\left(  \prod_{j=1}^{m-1}\left(
c_{j}d_{m}-c_{m}d_{j}\right)  \right)  \cdot\left(  \prod_{1\leq j<i\leq
m-1}\left(  c_{j}d_{i}-c_{i}d_{j}\right)  \right)  }_{\substack{=\prod_{1\leq
j<i\leq m}\left(  c_{j}d_{i}-c_{i}d_{j}\right)  \\\text{(by
(\ref{sol.cauchy-det.cauchy.a.5c}))}}}\\
&  =\underbrace{\left(  \dfrac{1}{x_{m,m}}\right)  ^{m-2}}_{=\dfrac{1}%
{x_{m,m}^{m-2}}}\underbrace{\left(  \dfrac{1}{x_{m,m}}\right)  \cdot x_{m,m}%
}_{=1}\cdot\dfrac{1}{\prod_{\left(  i,j\right)  \in\left\{  1,2,\ldots
,m\right\}  ^{2}}x_{i,j}}\\
&  \ \ \ \ \ \ \ \ \ \ \cdot\underbrace{\left(  \prod_{1\leq j<i\leq m}\left(
a_{i}b_{j}-a_{j}b_{i}\right)  \right)  \cdot\left(  \prod_{1\leq j<i\leq
m}\left(  c_{j}d_{i}-c_{i}d_{j}\right)  \right)  }_{=\prod_{1\leq j<i\leq
m}\left(  \left(  a_{i}b_{j}-a_{j}b_{i}\right)  \left(  c_{j}d_{i}-c_{i}%
d_{j}\right)  \right)  }\\
&  =\dfrac{1}{x_{m,m}^{m-2}}\cdot\dfrac{1}{\prod_{\left(  i,j\right)
\in\left\{  1,2,\ldots,m\right\}  ^{2}}x_{i,j}}\cdot\prod_{1\leq j<i\leq
m}\left(  \left(  a_{i}b_{j}-a_{j}b_{i}\right)  \left(  c_{j}d_{i}-c_{i}%
d_{j}\right)  \right)  .
\end{align*}
Multiplying both sides of this equality by $x_{m,m}^{m-2}$, we obtain%
\begin{align*}
&  x_{m,m}^{m-2}\cdot\det\left(  \left(  \dfrac{1}{x_{i,j}}\cdot\dfrac
{1}{x_{m,m}}-\dfrac{1}{x_{i,m}}\cdot\dfrac{1}{x_{m,j}}\right)  _{1\leq i\leq
m-1,\ 1\leq j\leq m-1}\right) \\
&  =\underbrace{x_{m,m}^{m-2}\cdot\dfrac{1}{x_{m,m}^{m-2}}}_{=1}\cdot\dfrac
{1}{\prod_{\left(  i,j\right)  \in\left\{  1,2,\ldots,m\right\}  ^{2}}x_{i,j}%
}\cdot\prod_{1\leq j<i\leq m}\left(  \left(  a_{i}b_{j}-a_{j}b_{i}\right)
\left(  c_{j}d_{i}-c_{i}d_{j}\right)  \right) \\
&  =\dfrac{1}{\prod_{\left(  i,j\right)  \in\left\{  1,2,\ldots,m\right\}
^{2}}x_{i,j}}\cdot\prod_{1\leq j<i\leq m}\left(  \left(  a_{i}b_{j}-a_{j}%
b_{i}\right)  \left(  c_{j}d_{i}-c_{i}d_{j}\right)  \right) \\
&  =\dfrac{\prod_{1\leq j<i\leq m}\left(  \left(  a_{i}b_{j}-a_{j}%
b_{i}\right)  \left(  c_{j}d_{i}-c_{i}d_{j}\right)  \right)  }{\prod_{\left(
i,j\right)  \in\left\{  1,2,\ldots,m\right\}  ^{2}}x_{i,j}}.
\end{align*}
Compared with (\ref{sol.cauchy-det.cauchy.a.1}), this yields
\[
\det\left(  \left(  \dfrac{1}{x_{i,j}}\right)  _{1\leq i\leq m,\ 1\leq j\leq
m}\right)  =\dfrac{\prod_{1\leq j<i\leq m}\left(  \left(  a_{i}b_{j}%
-a_{j}b_{i}\right)  \left(  c_{j}d_{i}-c_{i}d_{j}\right)  \right)  }%
{\prod_{\left(  i,j\right)  \in\left\{  1,2,\ldots,m\right\}  ^{2}}x_{i,j}}.
\]
In other words, (\ref{sol.cauchy-det.cauchy.a.goal}) holds for $k=m$. This
completes the induction step. The induction proof of
(\ref{sol.cauchy-det.cauchy.a.goal}) is thus complete.]

So we have proven (\ref{sol.cauchy-det.cauchy.a.goal}). We thus can apply
(\ref{sol.cauchy-det.cauchy.a.goal}) to $k=n$, and obtain%
\[
\det\left(  \left(  \dfrac{1}{x_{i,j}}\right)  _{1\leq i\leq n,\ 1\leq j\leq
n}\right)  =\dfrac{\prod_{1\leq j<i\leq n}\left(  \left(  a_{i}b_{j}%
-a_{j}b_{i}\right)  \left(  c_{j}d_{i}-c_{i}d_{j}\right)  \right)  }%
{\prod_{\left(  i,j\right)  \in\left\{  1,2,\ldots,n\right\}  ^{2}}x_{i,j}}.
\]
In light of (\ref{sol.cauchy-det.cauchy.xij=}), this rewrites as
\[
\det\left(  \left(  \dfrac{1}{a_{i}d_{j}-b_{i}c_{j}}\right)  _{1\leq i\leq
n,\ 1\leq j\leq n}\right)  =\dfrac{\prod_{1\leq j<i\leq n}\left(  \left(
a_{i}b_{j}-a_{j}b_{i}\right)  \left(  c_{j}d_{i}-c_{i}d_{j}\right)  \right)
}{\prod_{\left(  i,j\right)  \in\left\{  1,2,\ldots,n\right\}  ^{2}}\left(
a_{i}d_{j}-b_{i}c_{j}\right)  }.
\]
This proves Proposition \ref{prop.sol.cauchy-det.cauchy}.
\end{proof}

Now, let us solve Exercise \ref{exe.cauchy-det}.

\begin{proof}
[Solution to Exercise \ref{exe.cauchy-det}.]We have assumed that $x_{i}+y_{j}$
is an invertible element of $\mathbb{K}$ for every $\left(  i,j\right)
\in\left\{  1,2,\ldots,n\right\}  ^{2}$. In other words, $x_{i}+y_{j}$ is an
invertible element of $\mathbb{K}$ for every $i\in\left\{  1,2,\ldots
,n\right\}  $ and $j\in\left\{  1,2,\ldots,n\right\}  $. In other words,
$x_{i}\cdot1-\left(  -1\right)  \cdot y_{j}$ is an invertible element of
$\mathbb{K}$ for every $i\in\left\{  1,2,\ldots,n\right\}  $ and $j\in\left\{
1,2,\ldots,n\right\}  $ (because $\underbrace{x_{i}\cdot1}_{=x_{i}%
}-\underbrace{\left(  -1\right)  \cdot y_{j}}_{=-y_{j}}=x_{i}-\left(
-y_{j}\right)  =x_{i}+y_{j}$ for every $i\in\left\{  1,2,\ldots,n\right\}  $
and $j\in\left\{  1,2,\ldots,n\right\}  $). Hence, Proposition
\ref{prop.sol.cauchy-det.cauchy} (applied to $x_{i}$, $-1$, $y_{i}$ and $1$
instead of $a_{i}$, $b_{i}$, $c_{i}$ and $d_{i}$) yields%
\begin{align*}
&  \det\left(  \left(  \dfrac{1}{x_{i}\cdot1-\left(  -1\right)  \cdot y_{j}%
}\right)  _{1\leq i\leq n,\ 1\leq j\leq n}\right) \\
&  =\dfrac{\prod_{1\leq j<i\leq n}\left(  \left(  x_{i}\cdot\left(  -1\right)
-x_{j}\cdot\left(  -1\right)  \right)  \left(  y_{j}\cdot1-y_{i}\cdot1\right)
\right)  }{\prod_{\left(  i,j\right)  \in\left\{  1,2,\ldots,n\right\}  ^{2}%
}\left(  x_{i}\cdot1-\left(  -1\right)  \cdot y_{j}\right)  }\\
&  =\left(  \prod_{\left(  i,j\right)  \in\left\{  1,2,\ldots,n\right\}  ^{2}%
}\underbrace{\left(  x_{i}\cdot1-\left(  -1\right)  \cdot y_{j}\right)
}_{=x_{i}+y_{j}}\right)  ^{-1}\\
&  \ \ \ \ \ \ \ \ \ \ \cdot\prod_{1\leq j<i\leq n}\left(  \underbrace{\left(
x_{i}\cdot\left(  -1\right)  -x_{j}\cdot\left(  -1\right)  \right)  }%
_{=x_{j}-x_{i}}\underbrace{\left(  y_{j}\cdot1-y_{i}\cdot1\right)  }%
_{=y_{j}-y_{i}}\right) \\
&  =\left(  \prod_{\left(  i,j\right)  \in\left\{  1,2,\ldots,n\right\}  ^{2}%
}\left(  x_{i}+y_{j}\right)  \right)  ^{-1}\cdot\underbrace{\prod_{1\leq
j<i\leq n}\left(  \left(  x_{j}-x_{i}\right)  \left(  y_{j}-y_{i}\right)
\right)  }_{\substack{=\prod_{1\leq i<j\leq n}\left(  \left(  x_{i}%
-x_{j}\right)  \left(  y_{i}-y_{j}\right)  \right)  \\\text{(here, we have
renamed the index }\left(  j,i\right)  \\\text{as }\left(  i,j\right)  \text{
in the product)}}}\\
&  =\left(  \prod_{\left(  i,j\right)  \in\left\{  1,2,\ldots,n\right\}  ^{2}%
}\left(  x_{i}+y_{j}\right)  \right)  ^{-1}\cdot\prod_{1\leq i<j\leq n}\left(
\left(  x_{i}-x_{j}\right)  \left(  y_{i}-y_{j}\right)  \right) \\
&  =\frac{\prod_{1\leq i<j\leq n}\left(  \left(  x_{i}-x_{j}\right)  \left(
y_{i}-y_{j}\right)  \right)  }{\prod_{\left(  i,j\right)  \in\left\{
1,2,\ldots,n\right\}  ^{2}}\left(  x_{i}+y_{j}\right)  }.
\end{align*}
This rewrites as%
\[
\det\left(  \left(  \dfrac{1}{x_{i}+y_{j}}\right)  _{1\leq i\leq n,\ 1\leq
j\leq n}\right)  =\frac{\prod_{1\leq i<j\leq n}\left(  \left(  x_{i}%
-x_{j}\right)  \left(  y_{i}-y_{j}\right)  \right)  }{\prod_{\left(
i,j\right)  \in\left\{  1,2,\ldots,n\right\}  ^{2}}\left(  x_{i}+y_{j}\right)
}%
\]
(because every $i\in\left\{  1,2,\ldots,n\right\}  $ and $j\in\left\{
1,2,\ldots,n\right\}  $ satisfy $x_{i}\cdot1-\left(  -1\right)  \cdot
y_{j}=x_{i}+y_{j}$). This solves Exercise \ref{exe.cauchy-det}.
\end{proof}

\subsection{Solution to Exercise \ref{exe.det.schur-lem}}

\begin{proof}
[Solution to Exercise \ref{exe.det.schur-lem}.]For every $n\times n$-matrix
$\left(  c_{i,j}\right)  _{1\leq i\leq n,\ 1\leq j\leq n}$, we have%
\begin{equation}
\det\left(  \left(  c_{i,j}\right)  _{1\leq i\leq n,\ 1\leq j\leq n}\right)
=\sum\limits_{\sigma\in S_{n}}\left(  -1\right)  ^{\sigma}\prod\limits_{i=1}%
^{n}c_{\sigma\left(  i\right)  ,i} \label{sol.det.schur-lem.1}%
\end{equation}
\footnote{\textit{Proof of (\ref{sol.det.schur-lem.1}):} Let $\left(
c_{i,j}\right)  _{1\leq i\leq n,\ 1\leq j\leq n}$ be any $n\times n$-matrix.
Then, Exercise \ref{exe.ps4.4} (applied to $A=\left(  c_{i,j}\right)  _{1\leq
i\leq n,\ 1\leq j\leq n}$) yields%
\[
\det\left(  \left(  \left(  c_{i,j}\right)  _{1\leq i\leq n,\ 1\leq j\leq
n}\right)  ^{T}\right)  =\det\left(  \left(  c_{i,j}\right)  _{1\leq i\leq
n,\ 1\leq j\leq n}\right)  .
\]
Hence,%
\begin{align*}
\det\left(  \left(  c_{i,j}\right)  _{1\leq i\leq n,\ 1\leq j\leq n}\right)
&  =\det\left(  \underbrace{\left(  \left(  c_{i,j}\right)  _{1\leq i\leq
n,\ 1\leq j\leq n}\right)  ^{T}}_{\substack{=\left(  c_{j,i}\right)  _{1\leq
i\leq n,\ 1\leq j\leq n}\\\text{(by the definition of the}\\\text{transpose of
a matrix)}}}\right)  =\det\left(  \left(  c_{j,i}\right)  _{1\leq i\leq
n,\ 1\leq j\leq n}\right) \\
&  =\sum_{\sigma\in S_{n}}\left(  -1\right)  ^{\sigma}\prod_{i=1}^{n}%
c_{\sigma\left(  i\right)  ,i}\\
&  \ \ \ \ \ \ \ \ \ \ \left(
\begin{array}
[c]{c}%
\text{by (\ref{eq.det.eq.2}), applied to }\left(  c_{j,i}\right)  _{1\leq
i\leq n,\ 1\leq j\leq n}\text{ and }c_{j,i}\\
\text{instead of }A\text{ and }a_{i,j}%
\end{array}
\right)  .
\end{align*}
This proves (\ref{sol.det.schur-lem.1}).}. Applied to $\left(  c_{i,j}\right)
_{1\leq i\leq n,\ 1\leq j\leq n}=\left(  a_{i,j}\right)  _{1\leq i\leq
n,\ 1\leq j\leq n}$, this yields%
\begin{equation}
\det\left(  \left(  a_{i,j}\right)  _{1\leq i\leq n,\ 1\leq j\leq n}\right)
=\sum\limits_{\sigma\in S_{n}}\left(  -1\right)  ^{\sigma}\prod\limits_{i=1}%
^{n}a_{\sigma\left(  i\right)  ,i}. \label{sol.det.schur-lem.2}%
\end{equation}

Now, let $k\in\left\{  1,2,\ldots,n\right\}  $. For every $\sigma\in S_{n}$,
we have%
\begin{align}
\underbrace{\prod\limits_{i=1}^{n}}_{=\prod_{i\in\left\{  1,2,\ldots
,n\right\}  }}b_{\sigma\left(  i\right)  }^{\delta_{i,k}}  &  =\prod
_{i\in\left\{  1,2,\ldots,n\right\}  }b_{\sigma\left(  i\right)  }%
^{\delta_{i,k}}=\underbrace{b_{\sigma\left(  k\right)  }^{\delta_{k,k}}%
}_{\substack{=b_{\sigma\left(  k\right)  }^{1}\\\text{(since }\delta
_{k,k}=1\text{)}}}\prod\limits_{\substack{i\in\left\{  1,2,\ldots,n\right\}
;\\i\neq k}}\underbrace{b_{\sigma\left(  i\right)  }^{\delta_{i,k}}%
}_{\substack{=b_{\sigma\left(  i\right)  }^{0}\\\text{(since }\delta
_{i,k}=0\\\text{(since }i\neq k\text{))}}}\nonumber\\
&  \ \ \ \ \ \ \ \ \ \ \left(  \text{here, we have split off the factor for
}i=k\text{ from the product}\right) \nonumber\\
&  =\underbrace{b_{\sigma\left(  k\right)  }^{1}}_{=b_{\sigma\left(  k\right)
}}\prod\limits_{\substack{i\in\left\{  1,2,\ldots,n\right\}  ;\\i\neq
k}}\underbrace{b_{\sigma\left(  i\right)  }^{0}}_{=1}=b_{\sigma\left(
k\right)  }\underbrace{\prod\limits_{\substack{i\in\left\{  1,2,\ldots
,n\right\}  ;\\i\neq k}}1}_{=1}=b_{\sigma\left(  k\right)  }.
\label{sol.det.schur-lem.prod}%
\end{align}
But we can apply (\ref{sol.det.schur-lem.1}) to $\left(  c_{i,j}\right)
_{1\leq i\leq n,\ 1\leq j\leq n}=\left(  a_{i,j}b_{i}^{\delta_{j,k}}\right)
_{1\leq i\leq n,\ 1\leq j\leq n}$, and thus obtain%
\begin{align}
\det\left(  \left(  a_{i,j}b_{i}^{\delta_{j,k}}\right)  _{1\leq i\leq
n,\ 1\leq j\leq n}\right)   &  =\sum\limits_{\sigma\in S_{n}}\left(
-1\right)  ^{\sigma}\underbrace{\prod\limits_{i=1}^{n}\left(  a_{\sigma\left(
i\right)  ,i}b_{\sigma\left(  i\right)  }^{\delta_{i,k}}\right)  }_{=\left(
\prod\limits_{i=1}^{n}a_{\sigma\left(  i\right)  ,i}\right)  \left(
\prod\limits_{i=1}^{n}b_{\sigma\left(  i\right)  }^{\delta_{i,k}}\right)
}\nonumber\\
&  =\sum\limits_{\sigma\in S_{n}}\left(  -1\right)  ^{\sigma}\left(
\prod\limits_{i=1}^{n}a_{\sigma\left(  i\right)  ,i}\right)
\underbrace{\left(  \prod\limits_{i=1}^{n}b_{\sigma\left(  i\right)  }%
^{\delta_{i,k}}\right)  }_{\substack{=b_{\sigma\left(  k\right)  }\\\text{(by
(\ref{sol.det.schur-lem.prod}))}}}\nonumber\\
&  =\sum\limits_{\sigma\in S_{n}}\left(  -1\right)  ^{\sigma}\left(
\prod\limits_{i=1}^{n}a_{\sigma\left(  i\right)  ,i}\right)  b_{\sigma\left(
k\right)  }. \label{sol.det.schur-lem.k-th-term}%
\end{align}

Now, let us forget that we fixed $k$. We thus have proven
(\ref{sol.det.schur-lem.k-th-term}) for every $k\in\left\{  1,2,\ldots
,n\right\}  $. On the other hand, for every $\sigma\in S_{n}$, we have%
\begin{align}
\underbrace{\sum_{k=1}^{n}}_{=\sum_{k\in\left\{  1,2,\ldots,n\right\}  }%
}b_{\sigma\left(  k\right)  }  &  =\sum_{k\in\left\{  1,2,\ldots,n\right\}
}b_{\sigma\left(  k\right)  }=\underbrace{\sum_{k\in\left\{  1,2,\ldots
,n\right\}  }}_{=\sum_{k=1}^{n}}b_{k}\nonumber\\
&  \ \ \ \ \ \ \ \ \ \ \left(
\begin{array}
[c]{c}%
\text{here, we have substituted }k\text{ for }\sigma\left(  k\right)  \text{
in the sum,}\\
\text{since }\sigma\text{ is a bijection }\left\{  1,2,\ldots,n\right\}
\rightarrow\left\{  1,2,\ldots,n\right\} \\
\text{(since }\sigma\text{ is a permutation of }\left\{  1,2,\ldots,n\right\}
\text{ (since }\sigma\in S_{n}\text{))}%
\end{array}
\right) \nonumber\\
&  =\sum_{k=1}^{n}b_{k}=b_{1}+b_{2}+\cdots+b_{n}.
\label{sol.det.schur-lem.sum-over-k}%
\end{align}
Now,%
\begin{align*}
&  \sum\limits_{k=1}^{n}\underbrace{\det\left(  \left(  a_{i,j}b_{i}%
^{\delta_{j,k}}\right)  _{1\leq i\leq n,\ 1\leq j\leq n}\right)
}_{\substack{=\sum\limits_{\sigma\in S_{n}}\left(  -1\right)  ^{\sigma}\left(
\prod\limits_{i=1}^{n}a_{\sigma\left(  i\right)  ,i}\right)  b_{\sigma\left(
k\right)  }\\\text{(by (\ref{sol.det.schur-lem.k-th-term}))}}}\\
&  =\sum_{k=1}^{n}\sum\limits_{\sigma\in S_{n}}\left(  -1\right)  ^{\sigma
}\left(  \prod\limits_{i=1}^{n}a_{\sigma\left(  i\right)  ,i}\right)
b_{\sigma\left(  k\right)  }=\sum\limits_{\sigma\in S_{n}}\left(  -1\right)
^{\sigma}\left(  \prod\limits_{i=1}^{n}a_{\sigma\left(  i\right)  ,i}\right)
\underbrace{\sum_{k=1}^{n}b_{\sigma\left(  k\right)  }}_{\substack{=b_{1}%
+b_{2}+\cdots+b_{n}\\\text{(by (\ref{sol.det.schur-lem.sum-over-k}))}}}\\
&  =\sum\limits_{\sigma\in S_{n}}\left(  -1\right)  ^{\sigma}\left(
\prod\limits_{i=1}^{n}a_{\sigma\left(  i\right)  ,i}\right)  \left(
b_{1}+b_{2}+\cdots+b_{n}\right) \\
&  =\left(  b_{1}+b_{2}+\cdots+b_{n}\right)  \underbrace{\sum\limits_{\sigma
\in S_{n}}\left(  -1\right)  ^{\sigma}\prod\limits_{i=1}^{n}a_{\sigma\left(
i\right)  ,i}}_{\substack{=\det\left(  \left(  a_{i,j}\right)  _{1\leq i\leq
n,\ 1\leq j\leq n}\right)  \\\text{(by (\ref{sol.det.schur-lem.2}))}}}\\
&  =\left(  b_{1}+b_{2}+\cdots+b_{n}\right)  \det\left(  \left(
a_{i,j}\right)  _{1\leq i\leq n,\ 1\leq j\leq n}\right)  .
\end{align*}
This solves Exercise \ref{exe.det.schur-lem}.
\end{proof}

\subsection{\label{sect.sols.vander-det.s1.sol2}Second solution to Exercise
\ref{exe.vander-det.s1}}

Exercise \ref{exe.det.schur-lem} can be used to give a new (and simpler)
solution to Exercise \ref{exe.vander-det.s1}, suggested by Karthik Karnik:

\begin{vershort}
\begin{proof}
[Second solution to Exercise \ref{exe.vander-det.s1}.]Exercise
\ref{exe.det.schur-lem} (applied to $a_{i,j}=x_{i}^{n-j}$ and $b_{i}=x_{i}$)
shows that%
\begin{align}
&  \sum\limits_{k=1}^{n}\det\left(  \left(  x_{i}^{n-j}x_{i}^{\delta_{j,k}%
}\right)  _{1\leq i\leq n,\ 1\leq j\leq n}\right) \nonumber\\
&  =\left(  x_{1}+x_{2}+\cdots+x_{n}\right)  \underbrace{\det\left(  \left(
x_{i}^{n-j}\right)  _{1\leq i\leq n,\ 1\leq j\leq n}\right)  }%
_{\substack{=\prod_{1\leq i<j\leq n}\left(  x_{i}-x_{j}\right)  \\\text{(by
Theorem \ref{thm.vander-det} \textbf{(a)})}}}\nonumber\\
&  =\left(  x_{1}+x_{2}+\cdots+x_{n}\right)  \prod_{1\leq i<j\leq n}\left(
x_{i}-x_{j}\right)  . \label{sol.vander-det.s1.sol2.short.1}%
\end{align}

The left hand side of this equality is a sum of $n$ determinants. We shall now
show that $n-1$ of these determinants (namely, the ones that appear as addends
for $k>1$ in the sum) are $0$.

Indeed, every $k\in\left\{  2,3,\ldots,n\right\}  $ satisfies%
\begin{equation}
\det\left(  \left(  x_{i}^{n-j}x_{i}^{\delta_{j,k}}\right)  _{1\leq i\leq
n,\ 1\leq j\leq n}\right)  =0 \label{sol.vander-det.s1.sol2.short.2}%
\end{equation}
\footnote{\textit{Proof of (\ref{sol.vander-det.s1.sol2.short.2}):} Let
$k\in\left\{  2,3,\ldots,n\right\}  $. Let $A$ be the $n\times n$-matrix
$\left(  x_{i}^{n-j}x_{i}^{\delta_{j,k}}\right)  _{1\leq i\leq n,\ 1\leq j\leq
n}$.
\par
Let $h\in\left\{  1,2,\ldots,n\right\}  $. Since $A=\left(  x_{i}^{n-j}%
x_{i}^{\delta_{j,k}}\right)  _{1\leq i\leq n,\ 1\leq j\leq n}$, we have%
\begin{equation}
\left(  \text{the }\left(  h,k\right)  \text{-th entry of }A\right)
=x_{h}^{n-k}\underbrace{x_{h}^{\delta_{k,k}}}_{\substack{=x_{h}^{1}%
\\\text{(since }\delta_{k,k}=1\text{)}}}=x_{h}^{n-k}x_{h}^{1}=x_{h}^{n-k+1}.
\label{sol.vander-det.s1.sol2.short.2.pf.0}%
\end{equation}
On the other hand, $k-1\in\left\{  1,2,\ldots,n\right\}  $ (since
$k\in\left\{  2,3,\ldots,n\right\}  $), and thus the $\left(  h,k-1\right)
$-th entry of $A$ exists. Since $A=\left(  x_{i}^{n-j}x_{i}^{\delta_{j,k}%
}\right)  _{1\leq i\leq n,\ 1\leq j\leq n}$, this entry satisfies%
\[
\left(  \text{the }\left(  h,k-1\right)  \text{-th entry of }A\right)
=x_{h}^{n-\left(  k-1\right)  }\underbrace{x_{h}^{\delta_{k-1,k}}%
}_{\substack{=x_{h}^{0}\\\text{(since }\delta_{k-1,k}=0\text{)}}%
}=x_{h}^{n-\left(  k-1\right)  }x_{h}^{0}=x_{h}^{n-\left(  k-1\right)  }%
=x_{h}^{n-k+1}.
\]
Comparing this with (\ref{sol.vander-det.s1.sol2.short.2.pf.0}), we obtain%
\[
\left(  \text{the }\left(  h,k\right)  \text{-th entry of }A\right)  =\left(
\text{the }\left(  h,k-1\right)  \text{-th entry of }A\right)  .
\]
\par
Now, let us forget that we fixed $h$. We thus have shown that $\left(
\text{the }\left(  h,k\right)  \text{-th entry of }A\right)  =\left(
\text{the }\left(  h,k-1\right)  \text{-th entry of }A\right)  $ for every
$h\in\left\{  1,2,\ldots,n\right\}  $. In other words, the $k$-th column of
$A$ equals the $\left(  k-1\right)  $-st column of $A$. Thus, the matrix $A$
has two equal columns (since $k-1\neq k$). Therefore, Exercise \ref{exe.ps4.6}
\textbf{(f)} shows that $\det A=0$. Since $A=\left(  x_{i}^{n-j}x_{i}%
^{\delta_{j,k}}\right)  _{1\leq i\leq n,\ 1\leq j\leq n}$, this rewrites as
$\det\left(  \left(  x_{i}^{n-j}x_{i}^{\delta_{j,k}}\right)  _{1\leq i\leq
n,\ 1\leq j\leq n}\right)  =0$. Qed.}. Furthermore, every $\left(  i,j\right)
\in\left\{  1,2,\ldots,n\right\}  ^{2}$ satisfies%
\begin{equation}
x_{i}^{n-j}x_{i}^{\delta_{j,1}}=%
\begin{cases}
x_{i}^{n-j}, & \text{if }j>1;\\
x_{i}^{n}, & \text{if }j=1
\end{cases}
\label{sol.vander-det.s1.sol2.short.4}%
\end{equation}
\footnote{\textit{Proof of (\ref{sol.vander-det.s1.sol2.short.4}):} Let
$\left(  i,j\right)  \in\left\{  1,2,\ldots,n\right\}  ^{2}$. We need to prove
the equality (\ref{sol.vander-det.s1.sol2.short.4}). To do so, it clearly
suffices to prove the following two claims:
\par
\textit{Claim 1:} If $j>1$, then $x_{i}^{n-j}x_{i}^{\delta_{j,1}}=x_{i}^{n-j}%
$.
\par
\textit{Claim 2:} If $j=1$, then $x_{i}^{n-j}x_{i}^{\delta_{j,1}}=x_{i}^{n}$.
\par
\textit{Proof of Claim 1:} Assume that $j>1$. Thus, $j\neq1$, so that
$\delta_{j,1}=0$, so that $x_{i}^{\delta_{j,1}}=x_{i}^{0}=1$. Hence,
$x_{i}^{n-j}\underbrace{x_{i}^{\delta_{j,1}}}_{=1}=x_{i}^{n-j}$, and thus
Claim 1 is proven.
\par
\textit{Proof of Claim 2:} Assume that $j=1$. Thus, $\delta_{j,1}=1$, so that
$x_{i}^{\delta_{j,1}}=x_{i}^{1}$ and thus $x_{i}^{n-j}\underbrace{x_{i}%
^{\delta_{j,1}}}_{=x_{i}^{1}}=x_{i}^{n-j}x_{i}^{1}=x_{i}^{\left(  n-j\right)
+1}=x_{i}^{\left(  n-1\right)  +1}$ (since $j=1$), so that $x_{i}^{n-j}%
x_{i}^{\delta_{j,1}}=x_{i}^{\left(  n-1\right)  +1}=x_{i}^{n}$. This proves
Claim 2.
\par
Hence, (\ref{sol.vander-det.s1.sol2.short.4}) is proven (since we have proven
Claims 1 and 2).}. Now, (\ref{sol.vander-det.s1.sol2.short.1}) shows that%
\begin{align*}
&  \left(  x_{1}+x_{2}+\cdots+x_{n}\right)  \prod_{1\leq i<j\leq n}\left(
x_{i}-x_{j}\right) \\
&  =\sum\limits_{k=1}^{n}\det\left(  \left(  x_{i}^{n-j}x_{i}^{\delta_{j,k}%
}\right)  _{1\leq i\leq n,\ 1\leq j\leq n}\right) \\
&  =\det\left(  \left(  \underbrace{x_{i}^{n-j}x_{i}^{\delta_{j,1}}%
}_{\substack{=%
\begin{cases}
x_{i}^{n-j}, & \text{if }j>1;\\
x_{i}^{n}, & \text{if }j=1
\end{cases}
\\\text{(by (\ref{sol.vander-det.s1.sol2.short.4}))}}}\right)  _{1\leq i\leq
n,\ 1\leq j\leq n}\right)  +\sum\limits_{k=2}^{n}\underbrace{\det\left(
\left(  x_{i}^{n-j}x_{i}^{\delta_{j,k}}\right)  _{1\leq i\leq n,\ 1\leq j\leq
n}\right)  }_{\substack{=0\\\text{(by (\ref{sol.vander-det.s1.sol2.short.2}%
))}}}\\
&  \ \ \ \ \ \ \ \ \ \ \left(  \text{here, we have split off the addend for
}k=1\text{ from the sum}\right) \\
&  =\det\left(  \left(
\begin{cases}
x_{i}^{n-j}, & \text{if }j>1;\\
x_{i}^{n}, & \text{if }j=1
\end{cases}
\right)  _{1\leq i\leq n,\ 1\leq j\leq n}\right)  +\underbrace{\sum
\limits_{k=2}^{n}0}_{=0}\\
&  =\det\left(  \left(
\begin{cases}
x_{i}^{n-j}, & \text{if }j>1;\\
x_{i}^{n}, & \text{if }j=1
\end{cases}
\right)  _{1\leq i\leq n,\ 1\leq j\leq n}\right)  .
\end{align*}
This solves Exercise \ref{exe.vander-det.s1} again.
\end{proof}
\end{vershort}

\begin{verlong}
\begin{proof}
[Second solution to Exercise \ref{exe.vander-det.s1}.]Exercise
\ref{exe.det.schur-lem} (applied to $a_{i,j}=x_{i}^{n-j}$ and $b_{i}=x_{i}$)
shows that%
\begin{align}
&  \sum\limits_{k=1}^{n}\det\left(  \left(  x_{i}^{n-j}x_{i}^{\delta_{j,k}%
}\right)  _{1\leq i\leq n,\ 1\leq j\leq n}\right) \nonumber\\
&  =\left(  x_{1}+x_{2}+\cdots+x_{n}\right)  \underbrace{\det\left(  \left(
x_{i}^{n-j}\right)  _{1\leq i\leq n,\ 1\leq j\leq n}\right)  }%
_{\substack{=\prod_{1\leq i<j\leq n}\left(  x_{i}-x_{j}\right)  \\\text{(by
Theorem \ref{thm.vander-det} \textbf{(a)})}}}\nonumber\\
&  =\left(  x_{1}+x_{2}+\cdots+x_{n}\right)  \prod_{1\leq i<j\leq n}\left(
x_{i}-x_{j}\right)  . \label{sol.vander-det.s1.sol2.1}%
\end{align}

Now, every $k\in\left\{  2,3,\ldots,n\right\}  $ satisfies%
\begin{equation}
\det\left(  \left(  x_{i}^{n-j}x_{i}^{\delta_{j,k}}\right)  _{1\leq i\leq
n,\ 1\leq j\leq n}\right)  =0 \label{sol.vander-det.s1.sol2.2}%
\end{equation}
\footnote{\textit{Proof of (\ref{sol.vander-det.s1.sol2.2}):} Let
$k\in\left\{  2,3,\ldots,n\right\}  $. Let $A$ be the $n\times n$-matrix
$\left(  x_{i}^{n-j}x_{i}^{\delta_{j,k}}\right)  _{1\leq i\leq n,\ 1\leq j\leq
n}$. Thus,%
\begin{align}
\left(  \text{the }k\text{-th column of }A\right)   &  =\left(
\begin{array}
[c]{c}%
x_{1}^{n-k}x_{1}^{\delta_{k,k}}\\
x_{2}^{n-k}x_{2}^{\delta_{k,k}}\\
\vdots\\
x_{n}^{n-k}x_{n}^{\delta_{k,k}}%
\end{array}
\right)  =\left(  x_{i}^{n-k}\underbrace{x_{i}^{\delta_{k,k}}}%
_{\substack{=x_{i}^{1}\\\text{(since }\delta_{k,k}=1\\\text{(since
}k=k\text{))}}}\right)  _{1\leq i\leq n,\ 1\leq j\leq1}\nonumber\\
&  =\left(  \underbrace{x_{i}^{n-k}x_{i}^{1}}_{=x_{i}^{\left(  n-k\right)
+1}}\right)  _{1\leq i\leq n,\ 1\leq j\leq1}=\left(  x_{i}^{\left(
n-k\right)  +1}\right)  _{1\leq i\leq n,\ 1\leq j\leq1}.
\label{sol.vander-det.s1.sol2.2.pf.1}%
\end{align}
On the other hand, we have $k\in\left\{  2,3,\ldots,n\right\}  $, so that
$k-1\in\left\{  1,2,\ldots,n-1\right\}  \subseteq\left\{  1,2,\ldots
,n\right\}  $. Hence, the $\left(  k-1\right)  $-th column of $A$ exists. From
$A=\left(  x_{i}^{n-j}x_{i}^{\delta_{j,k}}\right)  _{1\leq i\leq n,\ 1\leq
j\leq n}$, we obtain%
\begin{align*}
\left(  \text{the }\left(  k-1\right)  \text{-th column of }A\right)   &
=\left(
\begin{array}
[c]{c}%
x_{1}^{n-\left(  k-1\right)  }x_{1}^{\delta_{k-1,k}}\\
x_{2}^{n-\left(  k-1\right)  }x_{2}^{\delta_{k-1,k}}\\
\vdots\\
x_{n}^{n-\left(  k-1\right)  }x_{n}^{\delta_{k-1,k}}%
\end{array}
\right)  =\left(  x_{i}^{n-\left(  k-1\right)  }\underbrace{x_{i}%
^{\delta_{k-1,k}}}_{\substack{=x_{i}^{0}\\\text{(since }\delta_{k-1,k}%
=0\\\text{(since }k-1\neq k\text{))}}}\right)  _{1\leq i\leq n,\ 1\leq j\leq
1}\\
&  =\left(  \underbrace{x_{i}^{n-\left(  k-1\right)  }}_{=x_{i}^{\left(
n-k\right)  +1}}\underbrace{x_{i}^{0}}_{=1}\right)  _{1\leq i\leq n,\ 1\leq
j\leq1}=\left(  x_{i}^{\left(  n-k\right)  +1}\right)  _{1\leq i\leq n,\ 1\leq
j\leq1}\\
&  =\left(  \text{the }k\text{-th column of }A\right)
\ \ \ \ \ \ \ \ \ \ \left(  \text{by (\ref{sol.vander-det.s1.sol2.2.pf.1}%
)}\right)  .
\end{align*}
Thus, the matrix $A$ has two equal columns (since $k-1\neq k$). Therefore,
Exercise \ref{exe.ps4.6} \textbf{(f)} shows that $\det A=0$. Since $A=\left(
x_{i}^{n-j}x_{i}^{\delta_{j,k}}\right)  _{1\leq i\leq n,\ 1\leq j\leq n}$,
this rewrites as $\det\left(  \left(  x_{i}^{n-j}x_{i}^{\delta_{j,k}}\right)
_{1\leq i\leq n,\ 1\leq j\leq n}\right)  =0$. Qed.}. Furthermore, every
$\left(  i,j\right)  \in\left\{  1,2,\ldots,n\right\}  ^{2}$ satisfies%
\begin{equation}
x_{i}^{n-j}x_{i}^{\delta_{j,1}}=%
\begin{cases}
x_{i}^{n-j}, & \text{if }j>1;\\
x_{i}^{n}, & \text{if }j=1
\end{cases}
\label{sol.vander-det.s1.sol2.4}%
\end{equation}
\footnote{\textit{Proof of (\ref{sol.vander-det.s1.sol2.4}):} Let $\left(
i,j\right)  \in\left\{  1,2,\ldots,n\right\}  ^{2}$. Thus, $i\in\left\{
1,2,\ldots,n\right\}  $ and $j\in\left\{  1,2,\ldots,n\right\}  $. We need to
prove the equality (\ref{sol.vander-det.s1.sol2.4}).
\par
We distinguish between the following two cases:
\par
\textit{Case 1:} We have $j\neq1$.
\par
\textit{Case 2:} We have $j=1$.
\par
Let us first consider Case 1. In this case, we have $j\neq1$. Thus,
$\delta_{j,1}=0$, so that $x_{i}^{\delta_{j,1}}=x_{i}^{0}=1$. Since
$j\in\left\{  1,2,\ldots,n\right\}  $ and $j\neq1$, we have $j\in\left\{
1,2,\ldots,n\right\}  \setminus\left\{  1\right\}  =\left\{  2,3,\ldots
,n\right\}  $, so that $j>1$. Thus, $%
\begin{cases}
x_{i}^{n-j}, & \text{if }j>1;\\
x_{i}^{n}, & \text{if }j=1
\end{cases}
=x_{i}^{n-j}$. Compared with $x_{i}^{n-j}\underbrace{x_{i}^{\delta_{j,1}}%
}_{=1}=x_{i}^{n-j}$, this yields $x_{i}^{n-j}x_{i}^{\delta_{j,1}}=%
\begin{cases}
x_{i}^{n-j}, & \text{if }j>1;\\
x_{i}^{n}, & \text{if }j=1
\end{cases}
$. Thus, (\ref{sol.vander-det.s1.sol2.4}) is proven in Case 1.
\par
Let us now consider Case 2. In this case, we have $j=1$. Thus, $\delta
_{j,1}=1$, so that $x_{i}^{\delta_{j,1}}=x_{i}^{1}$ and thus $x_{i}%
^{n-j}\underbrace{x_{i}^{\delta_{j,1}}}_{=x_{i}^{1}}=x_{i}^{n-j}x_{i}%
^{1}=x_{i}^{\left(  n-j\right)  +1}=x_{i}^{\left(  n-1\right)  +1}$ (since
$j=1$), so that $x_{i}^{n-j}x_{i}^{\delta_{j,1}}=x_{i}^{\left(  n-1\right)
+1}=x_{i}^{n}$. Compared with $%
\begin{cases}
x_{i}^{n-j}, & \text{if }j>1;\\
x_{i}^{n}, & \text{if }j=1
\end{cases}
=x_{i}^{n}$ (since $j=1$), this yields $x_{i}^{n-j}x_{i}^{\delta_{j,1}}=%
\begin{cases}
x_{i}^{n-j}, & \text{if }j>1;\\
x_{i}^{n}, & \text{if }j=1
\end{cases}
$. Thus, (\ref{sol.vander-det.s1.sol2.4}) is proven in Case 2.
\par
Now, we have proven (\ref{sol.vander-det.s1.sol2.4}) in each of the two Cases
1 and 2. Since these two Cases cover all possibilities, we thus see that
(\ref{sol.vander-det.s1.sol2.4}) always holds. Qed.}. Now,
(\ref{sol.vander-det.s1.sol2.1}) shows that%
\begin{align*}
&  \left(  x_{1}+x_{2}+\cdots+x_{n}\right)  \prod_{1\leq i<j\leq n}\left(
x_{i}-x_{j}\right) \\
&  =\sum\limits_{k=1}^{n}\det\left(  \left(  x_{i}^{n-j}x_{i}^{\delta_{j,k}%
}\right)  _{1\leq i\leq n,\ 1\leq j\leq n}\right) \\
&  =\det\left(  \left(  \underbrace{x_{i}^{n-j}x_{i}^{\delta_{j,1}}%
}_{\substack{=%
\begin{cases}
x_{i}^{n-j}, & \text{if }j>1;\\
x_{i}^{n}, & \text{if }j=1
\end{cases}
\\\text{(by (\ref{sol.vander-det.s1.sol2.4}))}}}\right)  _{1\leq i\leq
n,\ 1\leq j\leq n}\right)  +\sum\limits_{k=2}^{n}\underbrace{\det\left(
\left(  x_{i}^{n-j}x_{i}^{\delta_{j,k}}\right)  _{1\leq i\leq n,\ 1\leq j\leq
n}\right)  }_{\substack{=0\\\text{(by (\ref{sol.vander-det.s1.sol2.2}))}}}\\
&  \ \ \ \ \ \ \ \ \ \ \left(  \text{here, we have split off the addend for
}k=1\text{ from the sum}\right) \\
&  =\det\left(  \left(
\begin{cases}
x_{i}^{n-j}, & \text{if }j>1;\\
x_{i}^{n}, & \text{if }j=1
\end{cases}
\right)  _{1\leq i\leq n,\ 1\leq j\leq n}\right)  +\underbrace{\sum
\limits_{k=2}^{n}0}_{=0}\\
&  =\det\left(  \left(
\begin{cases}
x_{i}^{n-j}, & \text{if }j>1;\\
x_{i}^{n}, & \text{if }j=1
\end{cases}
\right)  _{1\leq i\leq n,\ 1\leq j\leq n}\right)  .
\end{align*}
This solves Exercise \ref{exe.vander-det.s1} again.
\end{proof}
\end{verlong}

\subsection{Solution to Exercise \ref{exe.det.a1a2anx}}

\begin{proof}
[First solution to Exercise \ref{exe.det.a1a2anx}.]The following solution will
rely on Exercise \ref{exe.ps4.6k} and \ref{exe.ps4.3}. Since this is not the
first time (nor the second) that we are using these exercises, I shall be brief.

Let $A$ be the $\left(  n+1\right)  \times\left(  n+1\right)  $-matrix
$\left(
\begin{array}
[c]{cccccc}%
x & a_{1} & a_{2} & \cdots & a_{n-1} & a_{n}\\
a_{1} & x & a_{2} & \cdots & a_{n-1} & a_{n}\\
a_{1} & a_{2} & x & \cdots & a_{n-1} & a_{n}\\
\vdots & \vdots & \vdots & \ddots & \vdots & \vdots\\
a_{1} & a_{2} & a_{3} & \cdots & x & a_{n}\\
a_{1} & a_{2} & a_{3} & \cdots & a_{n} & x
\end{array}
\right)  $. We need to prove that $\det A=\left(  x+\sum_{i=1}^{n}%
a_{i}\right)  \prod_{i=1}^{n}\left(  x-a_{i}\right)  $.

We perform the following operations on the matrix $A$ (in this order):

\begin{itemize}
\item We subtract the $2$-nd row from the $1$-st row.

\item We subtract the $3$-rd row from the $2$-nd row.

\item $\ldots$

\item We subtract the $\left(  n+1\right)  $-th row from the $n$-th row.
\end{itemize}

As we know from Exercise \ref{exe.ps4.6k} \textbf{(a)}, these operations do
not change the determinant of the matrix. Thus, if we denote by $B$ the result
of these operations, then $\det B=\det A$. On the other hand, it is easy to
see that%
\[
B=\left(
\begin{array}
[c]{cccccc}%
x-a_{1} & a_{1}-x & 0 & \cdots & 0 & 0\\
0 & x-a_{2} & a_{2}-x & \cdots & 0 & 0\\
0 & 0 & x-a_{3} & \cdots & 0 & 0\\
\vdots & \vdots & \vdots & \ddots & \vdots & \vdots\\
0 & 0 & 0 & \cdots & x-a_{n} & a_{n}-x\\
a_{1} & a_{2} & a_{3} & \cdots & a_{n} & x
\end{array}
\right)  .
\]
(Here, for every $i\in\left\{  1,2,\ldots,n\right\}  $, the $i$-th row of $B$
has $i$-th entry $x-a_{i}$ and $\left(  i+1\right)  $-th entry $a_{i}-x$; all
other entries in this row are $0$. The $\left(  n+1\right)  $-th row of $B$ is
$\left(  a_{1},a_{2},\ldots,a_{n},x\right)  $.)

Next, we apply the following operations to the matrix $B$ (in this order):

\begin{itemize}
\item We add the $1$-st column to the $2$-nd column.

\item We add the $2$-nd column to the $3$-rd column.

\item $\ldots$

\item We add the $n$-th column to the $\left(  n+1\right)  $-th column.
\end{itemize}

(Notice that the order in which we are performing these operations forces
their effects to accumulate; namely, at every step, the column that we are
adding has already been modified by the previous step.) As we know from
Exercise \ref{exe.ps4.6k} \textbf{(b)}, these operations do not change the
determinant of the matrix. Thus, if we denote by $C$ the result of these
operations, then $\det C=\det B=\det A$. On the other hand, it is easy to see
that%
\begin{align*}
&  C\\
&  =\left(
\begin{array}
[c]{cccccc}%
x-a_{1} & 0 & 0 & \cdots & 0 & 0\\
0 & x-a_{2} & 0 & \cdots & 0 & 0\\
0 & 0 & x-a_{3} & \cdots & 0 & 0\\
\vdots & \vdots & \vdots & \ddots & \vdots & \vdots\\
0 & 0 & 0 & \cdots & x-a_{n} & 0\\
a_{1} & a_{1}+a_{2} & a_{1}+a_{2}+a_{3} & \cdots & a_{1}+a_{2}+\cdots+a_{n} &
x+a_{1}+a_{2}+\cdots+a_{n}%
\end{array}
\right)  .
\end{align*}
This is a lower-triangular matrix. Thus, Exercise \ref{exe.ps4.3} shows that
its determinant is the product of its diagonal entries. In other words,%
\begin{align*}
\det C  &  =\left(  x-a_{1}\right)  \left(  x-a_{2}\right)  \cdots\left(
x-a_{n}\right)  \left(  x+a_{1}+a_{2}+\cdots+a_{n}\right) \\
&  =\left(  x+\sum_{i=1}^{n}a_{i}\right)  \prod_{i=1}^{n}\left(
x-a_{i}\right)  .
\end{align*}
Compared with $\det C=\det A$, this yields $\det A=\left(  x+\sum_{i=1}%
^{n}a_{i}\right)  \prod_{i=1}^{n}\left(  x-a_{i}\right)  $. This solves
Exercise \ref{exe.det.a1a2anx}.
\end{proof}

\begin{vershort}
\begin{proof}
[Second solution to Exercise \ref{exe.det.a1a2anx}.]For any $\left(
i,j\right)  \in\left\{  1,2,\ldots,n+1\right\}  ^{2}$, define an element
$u_{i,j}\in\mathbb{K}$ by%
\[
u_{i,j}=%
\begin{cases}
a_{j}, & \text{if }j<i;\\
x, & \text{if }j=i;\\
a_{j-1}, & \text{if }j>i
\end{cases}
.
\]
This $u_{i,j}$ is well-defined\footnote{\textit{Proof.} It is sufficient to
check that $a_{j}$ is well-defined when $j<i$, and that $a_{j-1}$ is
well-defined when $j>i$ (because $x$ is always well-defined). But this is
easy:
\par
\begin{itemize}
\item If $j<i$, then $j\in\left\{  1,2,\ldots,n\right\}  $ (since $j<i\leq
n+1$ yields $j\leq n$), and thus $a_{j}$ is well-defined.
\par
\item If $j>i$, then $j\in\left\{  2,3,\ldots,n+1\right\}  $ (since $j>i\geq1$
yields $j\geq2$), and thus $a_{j-1}$ is well-defined.
\end{itemize}
}. Now, we define an $\left(  n+1\right)  \times\left(  n+1\right)  $-matrix
$U$ by $U=\left(  u_{i,j}\right)  _{1\leq i\leq n+1,\ 1\leq j\leq n+1}$. Thus,%
\[
U=\left(  u_{i,j}\right)  _{1\leq i\leq n+1,\ 1\leq j\leq n+1}=\left(
\begin{array}
[c]{cccccc}%
x & a_{1} & a_{2} & \cdots & a_{n-1} & a_{n}\\
a_{1} & x & a_{2} & \cdots & a_{n-1} & a_{n}\\
a_{1} & a_{2} & x & \cdots & a_{n-1} & a_{n}\\
\vdots & \vdots & \vdots & \ddots & \vdots & \vdots\\
a_{1} & a_{2} & a_{3} & \cdots & x & a_{n}\\
a_{1} & a_{2} & a_{3} & \cdots & a_{n} & x
\end{array}
\right)  .
\]
Our goal is now to prove that $\det U=\left(  x+\sum_{i=1}^{n}a_{i}\right)
\prod_{i=1}^{n}\left(  x-a_{i}\right)  $.

For any $\left(  i,j\right)  \in\left\{  1,2,\ldots,n+1\right\}  ^{2}$, define
an element $s_{i,j}\in\mathbb{K}$ by%
\[
s_{i,j}=%
\begin{cases}
1, & \text{if }i\leq j;\\
0, & \text{if }i>j
\end{cases}
.
\]
Now, we define an $\left(  n+1\right)  \times\left(  n+1\right)  $-matrix $S$
by $S=\left(  s_{i,j}\right)  _{1\leq i\leq n+1,\ 1\leq j\leq n+1}$. The
matrix $S$ is upper-triangular\footnote{It looks as follows: $S=\left(
\begin{array}
[c]{cccc}%
1 & 1 & \cdots & 1\\
0 & 1 & \cdots & 1\\
\vdots & \vdots & \ddots & \vdots\\
0 & 0 & \cdots & 1
\end{array}
\right)  $.}. Since the determinant of an upper-triangular matrix is the
product of its diagonal entries, we thus obtain $\det S=1\cdot1\cdot
\cdots\cdot1=1$.

We extend the $n$-tuple $\left(  a_{1},a_{2},\ldots,a_{n}\right)
\in\mathbb{K}^{n}$ to an $\left(  n+1\right)  $-tuple $\left(  a_{1}%
,a_{2},\ldots,a_{n+1}\right)  \in\mathbb{K}^{n+1}$ by setting $a_{n+1}=0$.
Thus, an element $a_{k}\in\mathbb{K}$ is defined for every $k\in\left\{
1,2,\ldots,n+1\right\}  $.

Now, for every $\left(  i,j\right)  \in\left\{  1,2,\ldots,n+1\right\}  ^{2}$,
we have%
\begin{equation}
\sum_{k=1}^{n+1}u_{i,k}s_{k,j}=\sum_{k=1}^{j}a_{k}+s_{i,j}\left(
x-a_{j}\right)  \label{sol.det.a1a2anx.short.US}%
\end{equation}
\footnote{\textit{Proof of (\ref{sol.det.a1a2anx.short.US}):} Let $\left(
i,j\right)  \in\left\{  1,2,\ldots,n+1\right\}  ^{2}$. We need to prove the
equality (\ref{sol.det.a1a2anx.short.US}).
\par
We have $\left(  i,j\right)  \in\left\{  1,2,\ldots,n+1\right\}  ^{2}$. Thus,
$1\leq i\leq n+1$ and $1\leq j\leq n+1$. Now,%
\begin{align}
\sum_{k=1}^{n+1}u_{i,k}s_{k,j}  &  =\sum_{k=1}^{j}u_{i,k}\underbrace{s_{k,j}%
}_{\substack{=1\\\text{(by the definition}\\\text{of }s_{k,j}\text{, since
}k\leq j\text{)}}}+\sum_{k=j+1}^{n+1}u_{i,k}\underbrace{s_{k,j}}%
_{\substack{=0\\\text{(by the definition}\\\text{of }s_{k,j}\text{, since
}k>j\text{)}}}\ \ \ \ \ \ \ \ \ \ \left(  \text{since }1\leq j\leq n+1\right)
\nonumber\\
&  =\sum_{k=1}^{j}u_{i,k}1+\underbrace{\sum_{k=j+1}^{n+1}u_{i,k}0}_{=0}%
=\sum_{k=1}^{j}u_{i,k}1=\sum_{k=1}^{j}u_{i,k}.
\label{sol.det.a1a2anx.short.US.pf.1}%
\end{align}
\par
Now, we must be in one of the following two cases:
\par
\textit{Case 1:} We have $i\leq j$.
\par
\textit{Case 2:} We have $i>j$.
\par
Let us first consider Case 1. In this case, we have $i\leq j$. The definition
of $s_{i,j}$ shows that $s_{i,j}=1$ (since $i\leq j$). Therefore,%
\begin{equation}
\underbrace{\sum_{k=1}^{j}a_{k}}_{\substack{=\sum_{k=1}^{j-1}a_{k}%
+a_{j}\\\text{(here, we have split off the}\\\text{addend for }k=j\text{ from
the sum)}}}+\underbrace{s_{i,j}}_{=1}\left(  x-a_{j}\right)  =\sum_{k=1}%
^{j-1}a_{k}+a_{j}+\left(  x-a_{j}\right)  =\sum_{k=1}^{j-1}a_{k}+x.
\label{sol.det.a1a2anx.short.US.pf.3}%
\end{equation}
\par
We have $1\leq i\leq j$. Thus, $0\leq i-1\leq j-1$. Now,
(\ref{sol.det.a1a2anx.short.US.pf.1}) becomes%
\begin{align*}
\sum_{k=1}^{n+1}u_{i,k}s_{k,j}  &  =\sum_{k=1}^{j}u_{i,k}=\underbrace{\sum
_{k=1}^{i}u_{i,k}}_{\substack{=\sum_{k=1}^{i-1}u_{i,k}+u_{i,i}\\\text{(here,
we have split off the}\\\text{addend for }k=i\text{ from the sum)}}%
}+\sum_{k=i+1}^{j}u_{i,k}\ \ \ \ \ \ \ \ \ \ \left(  \text{since }1\leq i\leq
j\right) \\
&  =\sum_{k=1}^{i-1}\underbrace{u_{i,k}}_{\substack{=a_{k}\\\text{(by the
definition of }u_{i,k}\text{,}\\\text{since }k<i\text{)}}}+\underbrace{u_{i,i}%
}_{\substack{=x\\\text{(by the definition of }u_{i,i}\text{,}\\\text{since
}i=i\text{)}}}+\sum_{k=i+1}^{j}\underbrace{u_{i,k}}_{\substack{=a_{k-1}%
\\\text{(by the definition of }u_{i,k}\text{,}\\\text{since }k>i\text{)}}}\\
&  =\sum_{k=1}^{i-1}a_{k}+x+\underbrace{\sum_{k=i+1}^{j}a_{k-1}}%
_{\substack{=\sum_{k=i}^{j-1}a_{k}\\\text{(here, we have substituted }k\text{
for }k-1\text{ in the sum)}}}=\sum_{k=1}^{i-1}a_{k}+x+\sum_{k=i}^{j-1}a_{k}\\
&  =\underbrace{\sum_{k=1}^{i-1}a_{k}+\sum_{k=i}^{j-1}a_{k}}_{\substack{=\sum
_{k=1}^{j-1}a_{k}\\\text{(since }0\leq i-1\leq j-1\text{)}}}+x=\sum
_{k=1}^{j-1}a_{k}+x.
\end{align*}
Compared with (\ref{sol.det.a1a2anx.short.US.pf.3}), this yields%
\[
\sum_{k=1}^{n+1}u_{i,k}s_{k,j}=\sum_{k=1}^{j}a_{k}+s_{i,j}\left(
x-a_{j}\right)  .
\]
Thus, (\ref{sol.det.a1a2anx.short.US}) is proven in Case 1.
\par
Let us now consider Case 2. In this case, we have $i>j$. Hence, $j<i$. The
definition of $s_{i,j}$ therefore shows that $s_{i,j}=0$. Hence,%
\begin{equation}
\sum_{k=1}^{j}a_{k}+\underbrace{s_{i,j}}_{=0}\left(  x-a_{j}\right)
=\sum_{k=1}^{j}a_{k}+\underbrace{0\left(  x-a_{j}\right)  }_{=0}=\sum
_{k=1}^{j}a_{k}. \label{sol.det.a1a2anx.short.US.pf.5}%
\end{equation}
\par
But (\ref{sol.det.a1a2anx.short.US.pf.1}) becomes%
\[
\sum_{k=1}^{n+1}u_{i,k}s_{k,j}=\sum_{k=1}^{j}\underbrace{u_{i,k}%
}_{\substack{=a_{k}\\\text{(by the definition of }u_{i,k}\text{,}\\\text{since
}k\leq j<i\text{)}}}=\sum_{k=1}^{j}a_{k}.
\]
Compared with (\ref{sol.det.a1a2anx.short.US.pf.5}), this yields%
\[
\sum_{k=1}^{n+1}u_{i,k}s_{k,j}=\sum_{k=1}^{j}a_{k}+s_{i,j}\left(
x-a_{j}\right)  .
\]
Thus, (\ref{sol.det.a1a2anx.short.US}) is proven in Case 2.
\par
We have thus proven (\ref{sol.det.a1a2anx.short.US}) in each of the two Cases
1 and 2. Hence, (\ref{sol.det.a1a2anx.short.US}) always holds.}.

For any two objects $i$ and $j$, we define an element $\delta_{i,j}%
\in\mathbb{K}$ by $\delta_{i,j}=%
\begin{cases}
1, & \text{if }i=j;\\
0, & \text{if }i\neq j
\end{cases}
$. For any $\left(  i,j\right)  \in\left\{  1,2,\ldots,n+1\right\}  ^{2}$,
define an element $c_{i,j}\in\mathbb{K}$ by%
\[
c_{i,j}=\delta_{i,j}\left(  x-a_{j}\right)  +\delta_{i,n+1}\sum_{k=1}^{j}%
a_{k}.
\]
Let $C$ be the $\left(  n+1\right)  \times\left(  n+1\right)  $-matrix defined
by%
\[
C=\left(  c_{i,j}\right)  _{1\leq i\leq n+1,\ 1\leq j\leq n+1}.
\]

Now, for every $\left(  i,j\right)  \in\left\{  1,2,\ldots,n+1\right\}  ^{2}$,
we have%
\begin{equation}
\sum_{k=1}^{n+1}s_{i,k}c_{k,j}=\sum_{k=1}^{j}a_{k}+s_{i,j}\left(
x-a_{j}\right)  \label{sol.det.a1a2anx.short.SC}%
\end{equation}
\footnote{\textit{Proof of (\ref{sol.det.a1a2anx.short.SC}):} Let $\left(
i,j\right)  \in\left\{  1,2,\ldots,n+1\right\}  ^{2}$. We need to prove the
equality (\ref{sol.det.a1a2anx.short.SC}).
\par
We have $\left(  i,j\right)  \in\left\{  1,2,\ldots,n+1\right\}  ^{2}$. Thus,
$1\leq i\leq n+1$ and $1\leq j\leq n+1$.
\par
The definition of $s_{i,n+1}$ yields $s_{i,n+1}=1$ (since $i\leq n+1$). Now,%
\begin{align*}
&  \sum_{k=1}^{n+1}s_{i,k}c_{k,j}\\
&  =\underbrace{\sum_{r=1}^{n+1}}_{=\sum_{r\in\left\{  1,2,\ldots,n+1\right\}
}}s_{i,r}\underbrace{c_{r,j}}_{\substack{=\delta_{r,j}\left(  x-a_{j}\right)
+\delta_{r,n+1}\sum_{k=1}^{j}a_{k}\\\text{(by the definition of }%
c_{r,j}\text{)}}}\ \ \ \ \ \ \ \ \ \ \left(  \text{here, we renamed the
summation index }k\text{ as }r\right) \\
&  =\sum_{r\in\left\{  1,2,\ldots,n+1\right\}  }s_{i,r}\left(  \delta
_{r,j}\left(  x-a_{j}\right)  +\delta_{r,n+1}\sum_{k=1}^{j}a_{k}\right) \\
&  =\underbrace{\sum_{r\in\left\{  1,2,\ldots,n+1\right\}  }s_{i,r}%
\delta_{r,j}\left(  x-a_{j}\right)  }_{\substack{=s_{i,j}\delta_{j,j}\left(
x-a_{j}\right)  +\sum_{\substack{r\in\left\{  1,2,\ldots,n+1\right\}  ;\\r\neq
j}}s_{i,r}\delta_{r,j}\left(  x-a_{j}\right)  \\\text{(here, we have split off
the addend for }r=j\text{ from the sum)}}}+\underbrace{\sum_{r\in\left\{
1,2,\ldots,n+1\right\}  }s_{i,r}\delta_{r,n+1}\sum_{k=1}^{j}a_{k}%
}_{\substack{=s_{i,n+1}\delta_{n+1,n+1}\sum_{k=1}^{j}a_{k}+\sum
_{\substack{r\in\left\{  1,2,\ldots,n+1\right\}  ;\\r\neq n+1}}s_{i,r}%
\delta_{r,n+1}\sum_{k=1}^{j}a_{k}\\\text{(here, we have split off the addend
for }r=n+1\text{ from the sum)}}}\\
&  =s_{i,j}\underbrace{\delta_{j,j}}_{\substack{=1\\\text{(since }j=j\text{)}%
}}\left(  x-a_{j}\right)  +\sum_{\substack{r\in\left\{  1,2,\ldots
,n+1\right\}  ;\\r\neq j}}s_{i,r}\underbrace{\delta_{r,j}}%
_{\substack{=0\\\text{(since }r\neq j\text{)}}}\left(  x-a_{j}\right) \\
&  \ \ \ \ \ \ \ \ \ \ +\underbrace{s_{i,n+1}}_{=1}\underbrace{\delta
_{n+1,n+1}}_{\substack{=1\\\text{(since }n+1=n+1\text{)}}}\sum_{k=1}^{j}%
a_{k}+\sum_{\substack{r\in\left\{  1,2,\ldots,n+1\right\}  ;\\r\neq
n+1}}s_{i,r}\underbrace{\delta_{r,n+1}}_{\substack{=0\\\text{(since }r\neq
n+1\text{)}}}\sum_{k=1}^{j}a_{k}\\
&  =s_{i,j}\left(  x-a_{j}\right)  +\underbrace{\sum_{\substack{r\in\left\{
1,2,\ldots,n+1\right\}  ;\\r\neq j}}s_{i,r}0\left(  x-a_{j}\right)  }%
_{=0}+\sum_{k=1}^{j}a_{k}+\underbrace{\sum_{\substack{r\in\left\{
1,2,\ldots,n+1\right\}  ;\\r\neq n+1}}s_{i,r}0\sum_{k=1}^{j}a_{k}}_{=0}\\
&  =s_{i,j}\left(  x-a_{j}\right)  +\sum_{k=1}^{j}a_{k}=\sum_{k=1}^{j}%
a_{k}+s_{i,j}\left(  x-a_{j}\right)  .
\end{align*}
Thus, (\ref{sol.det.a1a2anx.short.SC}) is proven.}. Thus, for every $\left(
i,j\right)  \in\left\{  1,2,\ldots,n+1\right\}  ^{2}$, we have%
\begin{equation}
\sum_{k=1}^{n+1}s_{i,k}c_{k,j}=\sum_{k=1}^{j}a_{k}+s_{i,j}\left(
x-a_{j}\right)  =\sum_{k=1}^{n+1}u_{i,k}s_{k,j}
\label{sol.det.a1a2anx.short.SCvsUS}%
\end{equation}
(by (\ref{sol.det.a1a2anx.short.US})).

For every $i\in\left\{  1,2,\ldots,n\right\}  $, we have%
\begin{equation}
c_{i,i}=x-a_{i} \label{sol.det.a1a2anx.short.cii}%
\end{equation}
\footnote{\textit{Proof of (\ref{sol.det.a1a2anx.short.cii}):} Let
$i\in\left\{  1,2,\ldots,n\right\}  $. Thus, $i\neq n+1$, so that
$\delta_{i,n+1}=0$. Now, the definition of $c_{i,i}$ yields%
\[
c_{i,i}=\underbrace{\delta_{i,i}}_{\substack{=1\\\text{(since }i=i\text{)}%
}}\left(  x-a_{i}\right)  +\underbrace{\delta_{i,n+1}}_{=0}\sum_{k=1}^{i}%
a_{k}=\left(  x-a_{i}\right)  +\underbrace{0\sum_{k=1}^{i}a_{k}}_{=0}%
=x-a_{i}.
\]
This proves (\ref{sol.det.a1a2anx.short.cii}).}. Also,%
\begin{equation}
c_{n+1,n+1}=x+\sum_{i=1}^{n}a_{i} \label{sol.det.a1a2anx.short.cn+1n+1}%
\end{equation}
\footnote{\textit{Proof of (\ref{sol.det.a1a2anx.short.cn+1n+1}):} The
definition of $c_{n+1,n+1}$ yields%
\begin{align*}
c_{n+1,n+1}  &  =\underbrace{\delta_{n+1,n+1}}_{\substack{=1\\\text{(since
}n+1=n+1\text{)}}}\left(  x-a_{n+1}\right)  +\underbrace{\delta_{n+1,n+1}%
}_{\substack{=1\\\text{(since }n+1=n+1\text{)}}}\sum_{k=1}^{n+1}a_{k}=\left(
x-a_{n+1}\right)  +\sum_{k=1}^{n+1}a_{k}\\
&  =x+\underbrace{\left(  \sum_{k=1}^{n+1}a_{k}-a_{n+1}\right)  }_{=\sum
_{k=1}^{n}a_{k}}=x+\sum_{k=1}^{n}a_{k}=x+\sum_{i=1}^{n}a_{i}%
\end{align*}
(here, we have renamed the summation index $k$ as $i$). This proves
(\ref{sol.det.a1a2anx.short.cn+1n+1}).}.

But we have $c_{i,j}=0$ for every $\left(  i,j\right)  \in\left\{
1,2,\ldots,n+1\right\}  ^{2}$ satisfying $i<j$%
\ \ \ \ \footnote{\textit{Proof.} Let $\left(  i,j\right)  \in\left\{
1,2,\ldots,n+1\right\}  ^{2}$ be such that $i<j$. Thus, $1\leq i<j\leq n+1$,
so that $i<n+1$. Hence, $i\neq n+1$ and thus $\delta_{i,n+1}=0$. Also, $i<j$,
so that $i\neq j$ and thus $\delta_{i,j}=0$. Now, the definition of $c_{i,j}$
yields $c_{i,j}=\underbrace{\delta_{i,j}}_{=0}\left(  x-a_{j}\right)
+\underbrace{\delta_{i,n+1}}_{=0}\sum_{k=1}^{j}a_{k}=0\left(  x-a_{j}\right)
+0\sum_{k=1}^{j}a_{k}=0$, qed.}. Hence, Exercise \ref{exe.ps4.3} (applied to
$n+1$, $C$ and $c_{i,j}$ instead of $n$, $A$ and $a_{i,j}$) shows that%
\begin{align}
\det C  &  =c_{1,1}c_{2,2}\cdots c_{n+1,n+1}=\prod_{i=1}^{n+1}c_{i,i}=\left(
\prod_{i=1}^{n}\underbrace{c_{i,i}}_{\substack{=x-a_{i}\\\text{(by
(\ref{sol.det.a1a2anx.short.cii}))}}}\right)  \underbrace{c_{n+1,n+1}%
}_{\substack{=x+\sum_{i=1}^{n}a_{i}\\\text{(by
(\ref{sol.det.a1a2anx.short.cn+1n+1}))}}}\nonumber\\
&  \ \ \ \ \ \ \ \ \ \ \ \ \ \ \ \ \ \ \ \ \left(  \text{here, we have split
off the factor for }i=n+1\text{ from the product}\right) \nonumber\\
&  =\left(  \prod_{i=1}^{n}\left(  x-a_{i}\right)  \right)  \left(
x+\sum_{i=1}^{n}a_{i}\right) \nonumber\\
&  =\left(  x+\sum_{i=1}^{n}a_{i}\right)  \prod_{i=1}^{n}\left(
x-a_{i}\right)  . \label{sol.det.a1a2anx.short.detC}%
\end{align}

But recall that $S=\left(  s_{i,j}\right)  _{1\leq i\leq n+1,\ 1\leq j\leq
n+1}$ and $C=\left(  c_{i,j}\right)  _{1\leq i\leq n+1,\ 1\leq j\leq n+1}$.
Hence, the definition of the product of two matrices shows that
\begin{align}
SC  &  =\left(  \underbrace{\sum_{k=1}^{n+1}s_{i,k}c_{k,j}}_{\substack{=\sum
_{k=1}^{n+1}u_{i,k}s_{k,j}\\\text{(by (\ref{sol.det.a1a2anx.short.SCvsUS}))}%
}}\right)  _{1\leq i\leq n+1,\ 1\leq j\leq n+1}\nonumber\\
&  =\left(  \sum_{k=1}^{n+1}u_{i,k}s_{k,j}\right)  _{1\leq i\leq n+1,\ 1\leq
j\leq n+1}. \label{sol.det.a1a2anx.short.3}%
\end{align}

On the other hand, $U=\left(  u_{i,j}\right)  _{1\leq i\leq n+1,\ 1\leq j\leq
n+1}$ and $S=\left(  s_{i,j}\right)  _{1\leq i\leq n+1,\ 1\leq j\leq n+1}$.
Hence, the definition of the product of two matrices shows that
\[
US=\left(  \sum_{k=1}^{n+1}u_{i,k}s_{k,j}\right)  _{1\leq i\leq n+1,\ 1\leq
j\leq n+1}=SC\ \ \ \ \ \ \ \ \ \ \left(  \text{by
(\ref{sol.det.a1a2anx.short.3})}\right)  .
\]

Thus,%
\begin{align*}
\det\left(  \underbrace{US}_{=SC}\right)   &  =\det\left(  SC\right)
=\underbrace{\det S}_{=1}\cdot\det C\\
&  \ \ \ \ \ \ \ \ \ \ \left(
\begin{array}
[c]{c}%
\text{by Theorem \ref{thm.det(AB)}, applied to }n+1\text{, }S\text{ and }C\\
\text{instead of }n\text{, }A\text{ and }B
\end{array}
\right) \\
&  =\det C=\left(  x+\sum_{i=1}^{n}a_{i}\right)  \prod_{i=1}^{n}\left(
x-a_{i}\right)  \ \ \ \ \ \ \ \ \ \ \left(  \text{by
(\ref{sol.det.a1a2anx.short.detC})}\right)  .
\end{align*}
Compared with%
\begin{align*}
\det\left(  US\right)   &  =\det U\cdot\underbrace{\det S}_{=1}\\
&  \ \ \ \ \ \ \ \ \ \ \left(  \text{by Theorem \ref{thm.det(AB)}, applied to
}n+1\text{, }U\text{ and }S\text{ instead of }n\text{, }A\text{ and }B\right)
\\
&  =\det U,
\end{align*}
this yields
\[
\det U=\left(  x+\sum_{i=1}^{n}a_{i}\right)  \prod_{i=1}^{n}\left(
x-a_{i}\right)  .
\]
This solves Exercise \ref{exe.det.a1a2anx}.
\end{proof}
\end{vershort}

\begin{verlong}
\begin{proof}
[Second solution to Exercise \ref{exe.det.a1a2anx}.]For any $\left(
i,j\right)  \in\left\{  1,2,\ldots,n+1\right\}  ^{2}$, define an element
$u_{i,j}\in\mathbb{K}$ by%
\[
u_{i,j}=%
\begin{cases}
a_{j}, & \text{if }j<i;\\
x, & \text{if }j=i;\\
a_{j-1}, & \text{if }j>i
\end{cases}
.
\]
This $u_{i,j}$ is well-defined\footnote{\textit{Proof.} Let $\left(
i,j\right)  \in\left\{  1,2,\ldots,n+1\right\}  ^{2}$. Thus, $i\in\left\{
1,2,\ldots,n+1\right\}  $ and $j\in\left\{  1,2,\ldots,n+1\right\}  $. We must
prove that $u_{i,j}$ is well-defined.
\par
We are in one of the following three cases:
\par
\textit{Case 1:} We have $j<i$.
\par
\textit{Case 2:} We have $j=i$.
\par
\textit{Case 3:} We have $j>i$.
\par
Let us first consider Case 1. In this case, we have $j<i$. Thus, $j<i\leq n+1$
(since $i\in\left\{  1,2,\ldots,n+1\right\}  $), so that $j\leq\left(
n+1\right)  -1$ (since $j$ and $n+1$ are integers). Thus, $j\leq\left(
n+1\right)  -1=n$. Combined with $j\geq1$ (since $j\in\left\{  1,2,\ldots
,n+1\right\}  $), this shows that $1\leq j\leq n$. In other words,
$j\in\left\{  1,2,\ldots,n\right\}  $. Hence, $a_{j}$ is well-defined. Now,
the definition of $u_{i,j}$ says that $u_{i,j}=%
\begin{cases}
a_{j}, & \text{if }j<i;\\
x, & \text{if }j=i;\\
a_{j-1}, & \text{if }j>i
\end{cases}
=a_{j}$ (since $j<i$). Hence, $u_{i,j}$ is well-defined (since $a_{j}$ is
well-defined). Thus, in Case 1, we have proven that $u_{i,j}$ is well-defined.
\par
Let us next consider Case 2. In this case, we have $j=i$. Thus, the definition
of $u_{i,j}$ says that $u_{i,j}=%
\begin{cases}
a_{j}, & \text{if }j<i;\\
x, & \text{if }j=i;\\
a_{j-1}, & \text{if }j>i
\end{cases}
=x$ (since $j=i$). Hence, $u_{i,j}$ is well-defined (since $x$ is
well-defined). Thus, in Case 2, we have proven that $u_{i,j}$ is well-defined.
\par
Let us finally consider Case 3. In this case, we have $j>i$. Thus, $j>i\geq1$
(since $i\in\left\{  1,2,\ldots,n+1\right\}  $), so that $j-1>0$ and thus
$j-1\geq1$ (since $j-1$ is an integer). Also, $j\in\left\{  1,2,\ldots
,n+1\right\}  $, so that $j\leq n+1$ and thus $j-1\leq n$. Combined with
$j-1\geq1$, this yields $1\leq j-1\leq n$. In other words, $j-1\in\left\{
1,2,\ldots,n\right\}  $. Hence, $a_{j-1}$ is well-defined. Now, the definition
of $u_{i,j}$ says that $u_{i,j}=%
\begin{cases}
a_{j}, & \text{if }j<i;\\
x, & \text{if }j=i;\\
a_{j-1}, & \text{if }j>i
\end{cases}
=a_{j-1}$ (since $j>i$). Hence, $u_{i,j}$ is well-defined (since $a_{j-1}$ is
well-defined). Thus, in Case 3, we have proven that $u_{i,j}$ is well-defined.
\par
Now, in each of our three Cases 1, 2 and 3, we have shown that $u_{i,j}$ is
well-defined. Since these three Cases cover all possibilities, we thus
conclude that $u_{i,j}$ is always well-defined. Qed.}. Now, we define an
$\left(  n+1\right)  \times\left(  n+1\right)  $-matrix $U$ by $U=\left(
u_{i,j}\right)  _{1\leq i\leq n+1,\ 1\leq j\leq n+1}$. Thus,%
\begin{align*}
U  &  =\left(  \underbrace{u_{i,j}}_{=%
\begin{cases}
a_{j}, & \text{if }j<i;\\
x, & \text{if }j=i;\\
a_{j-1}, & \text{if }j>i
\end{cases}
}\right)  _{1\leq i\leq n+1,\ 1\leq j\leq n+1}=\left(
\begin{cases}
a_{j}, & \text{if }j<i;\\
x, & \text{if }j=i;\\
a_{j-1}, & \text{if }j>i
\end{cases}
\right)  _{1\leq i\leq n+1,\ 1\leq j\leq n+1}\\
&  =\left(
\begin{array}
[c]{cccccc}%
x & a_{1} & a_{2} & \cdots & a_{n-1} & a_{n}\\
a_{1} & x & a_{2} & \cdots & a_{n-1} & a_{n}\\
a_{1} & a_{2} & x & \cdots & a_{n-1} & a_{n}\\
\vdots & \vdots & \vdots & \ddots & \vdots & \vdots\\
a_{1} & a_{2} & a_{3} & \cdots & x & a_{n}\\
a_{1} & a_{2} & a_{3} & \cdots & a_{n} & x
\end{array}
\right)  .
\end{align*}
Our goal is now to prove that $\det U=\left(  x+\sum_{i=1}^{n}a_{i}\right)
\prod_{i=1}^{n}\left(  x-a_{i}\right)  $.

For any $\left(  i,j\right)  \in\left\{  1,2,\ldots,n+1\right\}  ^{2}$, define
an element $s_{i,j}\in\mathbb{K}$ by%
\[
s_{i,j}=%
\begin{cases}
1, & \text{if }i\leq j;\\
0, & \text{if }i>j
\end{cases}
.
\]
Now, we define an $\left(  n+1\right)  \times\left(  n+1\right)  $-matrix $S$
by $S=\left(  s_{i,j}\right)  _{1\leq i\leq n+1,\ 1\leq j\leq n+1}$. Then,
$\det S=1$\ \ \ \ \footnote{\textit{Proof.} We have $S=\left(  s_{i,j}\right)
_{1\leq i\leq n+1,\ 1\leq j\leq n+1}$. Therefore, the definition of the
transpose of a matrix yields $S^{T}=\left(  s_{j,i}\right)  _{1\leq i\leq
n+1,\ 1\leq j\leq n+1}$.
\par
Now, let $\left(  i,j\right)  \in\left\{  1,2,\ldots,n+1\right\}  ^{2}$ be
such that $i<j$. Then, $j>i$ (since $i<j$). Now, the definition of $s_{j,i}$
yields $s_{j,i}=%
\begin{cases}
1, & \text{if }j\leq i;\\
0, & \text{if }j>i
\end{cases}
=0$ (since $j>i$).
\par
Let us now forget that we fixed $\left(  i,j\right)  $. We thus have shown
that $s_{j,i}=0$ for every $\left(  i,j\right)  \in\left\{  1,2,\ldots
,n+1\right\}  ^{2}$ satisfying $i<j$. Therefore, Exercise \ref{exe.ps4.3}
(applied to $n+1$, $S^{T}$ and $s_{j,i}$ instead of $n$, $A$ and $a_{i,j}$)
shows that
\[
\det\left(  S^{T}\right)  =s_{1,1}s_{2,2}\cdots s_{n+1,n+1}=\prod_{i=1}%
^{n+1}\underbrace{s_{i,i}}_{\substack{=%
\begin{cases}
1, & \text{if }i\leq i;\\
0, & \text{if }i>i
\end{cases}
\\\text{(by the definition of }s_{i,i}\text{)}}}=\prod_{i=1}^{n+1}\underbrace{%
\begin{cases}
1, & \text{if }i\leq i;\\
0, & \text{if }i>i
\end{cases}
}_{\substack{=1\\\text{(since }i\leq i\text{)}}}=\prod_{i=1}^{n+1}1=1.
\]
But Exercise \ref{exe.ps4.4} (applied to $n+1$ and $S$ instead of $n$ and $A$)
shows that $\det\left(  S^{T}\right)  =\det S$. Comparing this with
$\det\left(  S^{T}\right)  =1$, we obtain $\det S=1$, qed.}.

We extend the $n$-tuple $\left(  a_{1},a_{2},\ldots,a_{n}\right)
\in\mathbb{K}^{n}$ to an $\left(  n+1\right)  $-tuple $\left(  a_{1}%
,a_{2},\ldots,a_{n+1}\right)  \in\mathbb{K}^{n+1}$ by setting $a_{n+1}=0$.
Thus, an element $a_{k}\in\mathbb{K}$ is defined for every $k\in\left\{
1,2,\ldots,n+1\right\}  $.

Now, for every $\left(  i,j\right)  \in\left\{  1,2,\ldots,n+1\right\}  ^{2}$,
we have%
\begin{equation}
\sum_{k=1}^{n+1}u_{i,k}s_{k,j}=\sum_{k=1}^{j}a_{k}+s_{i,j}\left(
x-a_{j}\right)  \label{sol.det.a1a2anx.US}%
\end{equation}
\footnote{\textit{Proof of (\ref{sol.det.a1a2anx.US}):} Let $\left(
i,j\right)  \in\left\{  1,2,\ldots,n+1\right\}  ^{2}$. We need to prove the
equality (\ref{sol.det.a1a2anx.US}).
\par
We have $\left(  i,j\right)  \in\left\{  1,2,\ldots,n+1\right\}  ^{2}$. Thus,
$i\in\left\{  1,2,\ldots,n+1\right\}  $ and $j\in\left\{  1,2,\ldots
,n+1\right\}  $. From $j\in\left\{  1,2,\ldots,n+1\right\}  $, we obtain
$1\leq j\leq n+1$. From $i\in\left\{  1,2,\ldots,n+1\right\}  $, we obtain
$1\leq i\leq n+1$. Now,%
\begin{align}
\sum_{k=1}^{n+1}u_{i,k}\underbrace{s_{k,j}}_{\substack{=%
\begin{cases}
1, & \text{if }k\leq j;\\
0, & \text{if }k>j
\end{cases}
\\\text{(by the definition of }s_{k,j}\text{)}}}  &  =\sum_{k=1}^{n+1}u_{i,k}%
\begin{cases}
1, & \text{if }k\leq j;\\
0, & \text{if }k>j
\end{cases}
\nonumber\\
&  =\sum_{k=1}^{j}u_{i,k}\underbrace{%
\begin{cases}
1, & \text{if }k\leq j;\\
0, & \text{if }k>j
\end{cases}
}_{\substack{=1\\\text{(since }k\leq j\text{)}}}+\sum_{k=j+1}^{n+1}%
u_{i,k}\underbrace{%
\begin{cases}
1, & \text{if }k\leq j;\\
0, & \text{if }k>j
\end{cases}
}_{\substack{=0\\\text{(since }k>j\text{)}}}\nonumber\\
&  \ \ \ \ \ \ \ \ \ \ \left(  \text{since }1\leq j\leq n+1\right) \nonumber\\
&  =\sum_{k=1}^{j}u_{i,k}1+\underbrace{\sum_{k=j+1}^{n+1}u_{i,k}0}_{=0}%
=\sum_{k=1}^{j}u_{i,k}1=\sum_{k=1}^{j}u_{i,k}. \label{sol.det.a1a2anx.US.pf.1}%
\end{align}
\par
Now, we must be in one of the following two cases:
\par
\textit{Case 1:} We have $i\leq j$.
\par
\textit{Case 2:} We have $i>j$.
\par
Let us first consider Case 1. In this case, we have $i\leq j$. The definition
of $s_{i,j}$ shows that $s_{i,j}=%
\begin{cases}
1, & \text{if }i\leq j;\\
0, & \text{if }i>j
\end{cases}
=1$ (since $i\leq j$). Therefore,%
\begin{equation}
\underbrace{\sum_{k=1}^{j}a_{k}}_{\substack{=\sum_{k=1}^{j-1}a_{k}%
+a_{j}\\\text{(here, we have split off the}\\\text{addend for }k=j\text{ from
the sum)}}}+\underbrace{s_{i,j}}_{=1}\left(  x-a_{j}\right)  =\sum_{k=1}%
^{j-1}a_{k}+a_{j}+\left(  x-a_{j}\right)  =\sum_{k=1}^{j-1}a_{k}+x.
\label{sol.det.a1a2anx.US.pf.3}%
\end{equation}
\par
We have $1\leq i\leq j$. Thus, $0\leq i-1\leq j-1$. Now,
(\ref{sol.det.a1a2anx.US.pf.1}) becomes%
\begin{align*}
\sum_{k=1}^{n+1}u_{i,k}s_{k,j}  &  =\sum_{k=1}^{j}u_{i,k}=\underbrace{\sum
_{k=1}^{i}u_{i,k}}_{\substack{=\sum_{k=1}^{i-1}u_{i,k}+u_{i,i}\\\text{(here,
we have split off the}\\\text{addend for }k=i\text{ from the sum)}}%
}+\sum_{k=i+1}^{j}u_{i,k}\ \ \ \ \ \ \ \ \ \ \left(  \text{since }1\leq i\leq
j\right) \\
&  =\sum_{k=1}^{i-1}\underbrace{u_{i,k}}_{\substack{=%
\begin{cases}
a_{k}, & \text{if }k<i;\\
x, & \text{if }k=i;\\
a_{k-1}, & \text{if }k>i
\end{cases}
\\\text{(by the definition of }u_{i,k}\text{)}}}+\underbrace{u_{i,i}%
}_{\substack{=%
\begin{cases}
a_{i}, & \text{if }i<i;\\
x, & \text{if }i=i;\\
a_{i-1}, & \text{if }i>i
\end{cases}
\\\text{(by the definition of }u_{i,i}\text{)}}}+\sum_{k=i+1}^{j}%
\underbrace{u_{i,k}}_{\substack{=%
\begin{cases}
a_{k}, & \text{if }k<i;\\
x, & \text{if }k=i;\\
a_{k-1}, & \text{if }k>i
\end{cases}
\\\text{(by the definition of }u_{i,k}\text{)}}}\\
&  =\sum_{k=1}^{i-1}\underbrace{%
\begin{cases}
a_{k}, & \text{if }k<i;\\
x, & \text{if }k=i;\\
a_{k-1}, & \text{if }k>i
\end{cases}
}_{\substack{=a_{k}\\\text{(since }k<i\text{)}}}+\underbrace{%
\begin{cases}
a_{i}, & \text{if }i<i;\\
x, & \text{if }i=i;\\
a_{i-1}, & \text{if }i>i
\end{cases}
}_{\substack{=x\\\text{(since }i=i\text{)}}}+\sum_{k=i+1}^{j}\underbrace{%
\begin{cases}
a_{k}, & \text{if }k<i;\\
x, & \text{if }k=i;\\
a_{k-1}, & \text{if }k>i
\end{cases}
}_{\substack{=a_{k-1}\\\text{(since }k>i\text{)}}}\\
&  =\sum_{k=1}^{i-1}a_{k}+x+\underbrace{\sum_{k=i+1}^{j}a_{k-1}}%
_{\substack{=\sum_{k=i}^{j-1}a_{k}\\\text{(here, we have substituted }k\text{
for }k-1\text{ in the sum)}}}=\sum_{k=1}^{i-1}a_{k}+x+\sum_{k=i}^{j-1}a_{k}\\
&  =\underbrace{\sum_{k=1}^{i-1}a_{k}+\sum_{k=i}^{j-1}a_{k}}_{\substack{=\sum
_{k=1}^{j-1}a_{k}\\\text{(since }0\leq i-1\leq j-1\text{)}}}+x=\sum
_{k=1}^{j-1}a_{k}+x.
\end{align*}
Compared with (\ref{sol.det.a1a2anx.US.pf.3}), this yields%
\[
\sum_{k=1}^{n+1}u_{i,k}s_{k,j}=\sum_{k=1}^{j}a_{k}+s_{i,j}\left(
x-a_{j}\right)  .
\]
Thus, (\ref{sol.det.a1a2anx.US}) is proven in Case 1.
\par
Let us now consider Case 2. In this case, we have $i>j$. Hence, $j<i$. The
definition of $s_{i,j}$ shows that $s_{i,j}=%
\begin{cases}
1, & \text{if }i\leq j;\\
0, & \text{if }i>j
\end{cases}
=0$ (since $i>j$). Therefore,%
\begin{equation}
\sum_{k=1}^{j}a_{k}+\underbrace{s_{i,j}}_{=0}\left(  x-a_{j}\right)
=\sum_{k=1}^{j}a_{k}+\underbrace{0\left(  x-a_{j}\right)  }_{=0}=\sum
_{k=1}^{j}a_{k}. \label{sol.det.a1a2anx.US.pf.5}%
\end{equation}
\par
But (\ref{sol.det.a1a2anx.US.pf.1}) becomes%
\[
\sum_{k=1}^{n+1}u_{i,k}s_{k,j}=\sum_{k=1}^{j}\underbrace{u_{i,k}}_{\substack{=%
\begin{cases}
a_{k}, & \text{if }k<i;\\
x, & \text{if }k=i;\\
a_{k-1}, & \text{if }k>i
\end{cases}
\\\text{(by the definition of }u_{i,k}\text{)}}}=\sum_{k=1}^{j}\underbrace{%
\begin{cases}
a_{k}, & \text{if }k<i;\\
x, & \text{if }k=i;\\
a_{k-1}, & \text{if }k>i
\end{cases}
}_{\substack{=a_{k}\\\text{(since }k\leq j<i\text{)}}}=\sum_{k=1}^{j}a_{k}.
\]
Compared with (\ref{sol.det.a1a2anx.US.pf.5}), this yields%
\[
\sum_{k=1}^{n+1}u_{i,k}s_{k,j}=\sum_{k=1}^{j}a_{k}+s_{i,j}\left(
x-a_{j}\right)  .
\]
Thus, (\ref{sol.det.a1a2anx.US}) is proven in Case 2.
\par
We have thus proven (\ref{sol.det.a1a2anx.US}) in each of the two Cases 1 and
2. Since these two Cases cover all possibilities, this shows that
(\ref{sol.det.a1a2anx.US}) always holds. Qed.}.

For any two objects $i$ and $j$, we define an element $\delta_{i,j}%
\in\mathbb{K}$ by $\delta_{i,j}=%
\begin{cases}
1, & \text{if }i=j;\\
0, & \text{if }i\neq j
\end{cases}
$. For any $\left(  i,j\right)  \in\left\{  1,2,\ldots,n+1\right\}  ^{2}$,
define an element $c_{i,j}\in\mathbb{K}$ by%
\[
c_{i,j}=\delta_{i,j}\left(  x-a_{j}\right)  +\delta_{i,n+1}\sum_{k=1}^{j}%
a_{k}.
\]
Let $C$ be the $\left(  n+1\right)  \times\left(  n+1\right)  $-matrix defined
by%
\[
C=\left(  c_{i,j}\right)  _{1\leq i\leq n+1,\ 1\leq j\leq n+1}.
\]

Now, for every $\left(  i,j\right)  \in\left\{  1,2,\ldots,n+1\right\}  ^{2}$,
we have%
\begin{equation}
\sum_{k=1}^{n+1}s_{i,k}c_{k,j}=\sum_{k=1}^{j}a_{k}+s_{i,j}\left(
x-a_{j}\right)  \label{sol.det.a1a2anx.SC}%
\end{equation}
\footnote{\textit{Proof of (\ref{sol.det.a1a2anx.SC}):} Let $\left(
i,j\right)  \in\left\{  1,2,\ldots,n+1\right\}  ^{2}$. We need to prove the
equality (\ref{sol.det.a1a2anx.SC}).
\par
We have $\left(  i,j\right)  \in\left\{  1,2,\ldots,n+1\right\}  ^{2}$. Thus,
$i\in\left\{  1,2,\ldots,n+1\right\}  $ and $j\in\left\{  1,2,\ldots
,n+1\right\}  $. From $j\in\left\{  1,2,\ldots,n+1\right\}  $, we obtain
$1\leq j\leq n+1$. From $i\in\left\{  1,2,\ldots,n+1\right\}  $, we obtain
$1\leq i\leq n+1$.
\par
The definition of $s_{i,n+1}$ yields $s_{i,n+1}=%
\begin{cases}
1, & \text{if }i\leq n+1;\\
0, & \text{if }i>n+1
\end{cases}
=1$ (since $i\leq n+1$). Now,%
\begin{align*}
&  \sum_{k=1}^{n+1}s_{i,k}c_{k,j}\\
&  =\underbrace{\sum_{r=1}^{n+1}}_{=\sum_{r\in\left\{  1,2,\ldots,n+1\right\}
}}s_{i,r}\underbrace{c_{r,j}}_{\substack{=\delta_{r,j}\left(  x-a_{j}\right)
+\delta_{r,n+1}\sum_{k=1}^{j}a_{k}\\\text{(by the definition of }%
c_{r,j}\text{)}}}\ \ \ \ \ \ \ \ \ \ \left(  \text{here, we renamed the
summation index }k\text{ as }r\right) \\
&  =\sum_{r\in\left\{  1,2,\ldots,n+1\right\}  }s_{i,r}\left(  \delta
_{r,j}\left(  x-a_{j}\right)  +\delta_{r,n+1}\sum_{k=1}^{j}a_{k}\right) \\
&  =\underbrace{\sum_{r\in\left\{  1,2,\ldots,n+1\right\}  }s_{i,r}%
\delta_{r,j}\left(  x-a_{j}\right)  }_{\substack{=s_{i,j}\delta_{j,j}\left(
x-a_{j}\right)  +\sum_{\substack{r\in\left\{  1,2,\ldots,n+1\right\}  ;\\r\neq
j}}s_{i,r}\delta_{r,j}\left(  x-a_{j}\right)  \\\text{(here, we have split off
the addend for }r=j\text{ from the sum)}}}+\underbrace{\sum_{r\in\left\{
1,2,\ldots,n+1\right\}  }s_{i,r}\delta_{r,n+1}\sum_{k=1}^{j}a_{k}%
}_{\substack{=s_{i,n+1}\delta_{n+1,n+1}\sum_{k=1}^{j}a_{k}+\sum
_{\substack{r\in\left\{  1,2,\ldots,n+1\right\}  ;\\r\neq n+1}}s_{i,r}%
\delta_{r,n+1}\sum_{k=1}^{j}a_{k}\\\text{(here, we have split off the addend
for }r=n+1\text{ from the sum)}}}\\
&  =s_{i,j}\underbrace{\delta_{j,j}}_{\substack{=1\\\text{(since }j=j\text{)}%
}}\left(  x-a_{j}\right)  +\sum_{\substack{r\in\left\{  1,2,\ldots
,n+1\right\}  ;\\r\neq j}}s_{i,r}\underbrace{\delta_{r,j}}%
_{\substack{=0\\\text{(since }r\neq j\text{)}}}\left(  x-a_{j}\right) \\
&  \ \ \ \ \ \ \ \ \ \ +\underbrace{s_{i,n+1}}_{=1}\underbrace{\delta
_{n+1,n+1}}_{\substack{=1\\\text{(since }n+1=n+1\text{)}}}\sum_{k=1}^{j}%
a_{k}+\sum_{\substack{r\in\left\{  1,2,\ldots,n+1\right\}  ;\\r\neq
n+1}}s_{i,r}\underbrace{\delta_{r,n+1}}_{\substack{=0\\\text{(since }r\neq
n+1\text{)}}}\sum_{k=1}^{j}a_{k}\\
&  =s_{i,j}\left(  x-a_{j}\right)  +\underbrace{\sum_{\substack{r\in\left\{
1,2,\ldots,n+1\right\}  ;\\r\neq j}}s_{i,r}0\left(  x-a_{j}\right)  }%
_{=0}+\sum_{k=1}^{j}a_{k}+\underbrace{\sum_{\substack{r\in\left\{
1,2,\ldots,n+1\right\}  ;\\r\neq n+1}}s_{i,r}0\sum_{k=1}^{j}a_{k}}_{=0}\\
&  =s_{i,j}\left(  x-a_{j}\right)  +\sum_{k=1}^{j}a_{k}=\sum_{k=1}^{j}%
a_{k}+s_{i,j}\left(  x-a_{j}\right)  .
\end{align*}
Thus, (\ref{sol.det.a1a2anx.SC}) is proven.}. Thus, for every $\left(
i,j\right)  \in\left\{  1,2,\ldots,n+1\right\}  ^{2}$, we have%
\begin{equation}
\sum_{k=1}^{n+1}s_{i,k}c_{k,j}=\sum_{k=1}^{j}a_{k}+s_{i,j}\left(
x-a_{j}\right)  =\sum_{k=1}^{n+1}u_{i,k}s_{k,j} \label{sol.det.a1a2anx.SCvsUS}%
\end{equation}
(by (\ref{sol.det.a1a2anx.US})).

For every $i\in\left\{  1,2,\ldots,n\right\}  $, we have%
\begin{equation}
c_{i,i}=x-a_{i} \label{sol.det.a1a2anx.cii}%
\end{equation}
\footnote{\textit{Proof of (\ref{sol.det.a1a2anx.cii}):} Let $i\in\left\{
1,2,\ldots,n\right\}  $. Thus, $1\leq i\leq n$, so that $i\leq n<n+1$, and
therefore $i\neq n+1$. Hence, $\delta_{i,n+1}=0$. Now, the definition of
$c_{i,i}$ yields%
\[
c_{i,i}=\underbrace{\delta_{i,i}}_{\substack{=1\\\text{(since }i=i\text{)}%
}}\left(  x-a_{i}\right)  +\underbrace{\delta_{i,n+1}}_{=0}\sum_{k=1}^{i}%
a_{k}=\left(  x-a_{i}\right)  +\underbrace{0\sum_{k=1}^{i}a_{k}}_{=0}%
=x-a_{i}.
\]
This proves (\ref{sol.det.a1a2anx.cii}).}. Also,%
\begin{equation}
c_{n+1,n+1}=x+\sum_{i=1}^{n}a_{i} \label{sol.det.a1a2anx.cn+1n+1}%
\end{equation}
\footnote{\textit{Proof of (\ref{sol.det.a1a2anx.cn+1n+1}):} The definition of
$c_{n+1,n+1}$ yields%
\begin{align*}
c_{n+1,n+1}  &  =\underbrace{\delta_{n+1,n+1}}_{\substack{=1\\\text{(since
}n+1=n+1\text{)}}}\left(  x-a_{n+1}\right)  +\underbrace{\delta_{n+1,n+1}%
}_{\substack{=1\\\text{(since }n+1=n+1\text{)}}}\sum_{k=1}^{n+1}a_{k}=\left(
x-a_{n+1}\right)  +\underbrace{\sum_{k=1}^{n+1}a_{k}}_{\substack{=\sum
_{k=1}^{n}a_{k}+a_{n+1}\\\text{(here, we have split off the}\\\text{addend for
}k=n+1\text{ from the sum)}}}\\
&  =\left(  x-a_{n+1}\right)  +\sum_{k=1}^{n}a_{k}+a_{n+1}=x+\sum_{k=1}%
^{n}a_{k}=x+\sum_{i=1}^{n}a_{i}%
\end{align*}
(here, we have renamed the summation index $k$ as $i$). This proves
(\ref{sol.det.a1a2anx.cn+1n+1}).}.

But we have $c_{i,j}=0$ for every $\left(  i,j\right)  \in\left\{
1,2,\ldots,n+1\right\}  ^{2}$ satisfying $i<j$%
\ \ \ \ \footnote{\textit{Proof.} Let $\left(  i,j\right)  \in\left\{
1,2,\ldots,n+1\right\}  ^{2}$ be such that $i<j$. We have $\left(  i,j\right)
\in\left\{  1,2,\ldots,n+1\right\}  ^{2}$, so that $i\in\left\{
1,2,\ldots,n+1\right\}  $ and $j\in\left\{  1,2,\ldots,n+1\right\}  $. From
$j\in\left\{  1,2,\ldots,n+1\right\}  $, we obtain $1\leq j\leq n+1$. Thus,
$i<j\leq n+1$, so that $i\neq n+1$ and thus $\delta_{i,n+1}=0$. Also, $i<j$,
so that $i\neq j$ and thus $\delta_{i,j}=0$. Now, the definition of $c_{i,j}$
yields $c_{i,j}=\underbrace{\delta_{i,j}}_{=0}\left(  x-a_{j}\right)
+\underbrace{\delta_{i,n+1}}_{=0}\sum_{k=1}^{j}a_{k}=0\left(  x-a_{j}\right)
+0\sum_{k=1}^{j}a_{k}=0$, qed.}. Hence, Exercise \ref{exe.ps4.3} (applied to
$n+1$, $C$ and $c_{i,j}$ instead of $n$, $A$ and $a_{i,j}$) shows that%
\begin{align}
\det C  &  =c_{1,1}c_{2,2}\cdots c_{n+1,n+1}=\prod_{i=1}^{n+1}c_{i,i}%
\nonumber\\
&  =\left(  \prod_{i=1}^{n}\underbrace{c_{i,i}}_{\substack{=x-a_{i}\\\text{(by
(\ref{sol.det.a1a2anx.cii}))}}}\right)  \underbrace{c_{n+1,n+1}}%
_{\substack{=x+\sum_{i=1}^{n}a_{i}\\\text{(by (\ref{sol.det.a1a2anx.cn+1n+1}%
))}}}\ \ \ \ \ \ \ \ \ \ \left(
\begin{array}
[c]{c}%
\text{here, we have split off the factor}\\
\text{for }i=n+1\text{ from the product}%
\end{array}
\right) \nonumber\\
&  =\left(  \prod_{i=1}^{n}\left(  x-a_{i}\right)  \right)  \left(
x+\sum_{i=1}^{n}a_{i}\right) \nonumber\\
&  =\left(  x+\sum_{i=1}^{n}a_{i}\right)  \prod_{i=1}^{n}\left(
x-a_{i}\right)  . \label{sol.det.a1a2anx.detC}%
\end{align}

But recall that $S=\left(  s_{i,j}\right)  _{1\leq i\leq n+1,\ 1\leq j\leq
n+1}$ and $C=\left(  c_{i,j}\right)  _{1\leq i\leq n+1,\ 1\leq j\leq n+1}$.
Hence, the definition of the product of two matrices shows that
\begin{align}
SC  &  =\left(  \underbrace{\sum_{k=1}^{n+1}s_{i,k}c_{k,j}}_{\substack{=\sum
_{k=1}^{n+1}u_{i,k}s_{k,j}\\\text{(by (\ref{sol.det.a1a2anx.SCvsUS}))}%
}}\right)  _{1\leq i\leq n+1,\ 1\leq j\leq n+1}\nonumber\\
&  =\left(  \sum_{k=1}^{n+1}u_{i,k}s_{k,j}\right)  _{1\leq i\leq n+1,\ 1\leq
j\leq n+1}. \label{sol.det.a1a2anx.3}%
\end{align}

On the other hand, $U=\left(  u_{i,j}\right)  _{1\leq i\leq n+1,\ 1\leq j\leq
n+1}$ and $S=\left(  s_{i,j}\right)  _{1\leq i\leq n+1,\ 1\leq j\leq n+1}$.
Hence, the definition of the product of two matrices shows that
\[
US=\left(  \sum_{k=1}^{n+1}u_{i,k}s_{k,j}\right)  _{1\leq i\leq n+1,\ 1\leq
j\leq n+1}=SC\ \ \ \ \ \ \ \ \ \ \left(  \text{by (\ref{sol.det.a1a2anx.3}%
)}\right)  .
\]

Thus,%
\begin{align*}
\det\left(  \underbrace{US}_{=SC}\right)   &  =\det\left(  SC\right)
=\underbrace{\det S}_{=1}\cdot\det C\\
&  \ \ \ \ \ \ \ \ \ \ \left(  \text{by Theorem \ref{thm.det(AB)}, applied to
}n+1\text{, }S\text{ and }C\text{ instead of }n\text{, }A\text{ and }B\right)
\\
&  =\det C=\left(  x+\sum_{i=1}^{n}a_{i}\right)  \prod_{i=1}^{n}\left(
x-a_{i}\right)  \ \ \ \ \ \ \ \ \ \ \left(  \text{by
(\ref{sol.det.a1a2anx.detC})}\right)  .
\end{align*}
Compared with%
\begin{align*}
\det\left(  US\right)   &  =\det U\cdot\underbrace{\det S}_{=1}\\
&  \ \ \ \ \ \ \ \ \ \ \left(  \text{by Theorem \ref{thm.det(AB)}, applied to
}n+1\text{, }U\text{ and }S\text{ instead of }n\text{, }A\text{ and }B\right)
\\
&  =\det U,
\end{align*}
this yields
\[
\det U=\left(  x+\sum_{i=1}^{n}a_{i}\right)  \prod_{i=1}^{n}\left(
x-a_{i}\right)  .
\]
Since $U=\left(
\begin{array}
[c]{cccccc}%
x & a_{1} & a_{2} & \cdots & a_{n-1} & a_{n}\\
a_{1} & x & a_{2} & \cdots & a_{n-1} & a_{n}\\
a_{1} & a_{2} & x & \cdots & a_{n-1} & a_{n}\\
\vdots & \vdots & \vdots & \ddots & \vdots & \vdots\\
a_{1} & a_{2} & a_{3} & \cdots & x & a_{n}\\
a_{1} & a_{2} & a_{3} & \cdots & a_{n} & x
\end{array}
\right)  $, this rewrites as%
\[
\det\left(
\begin{array}
[c]{cccccc}%
x & a_{1} & a_{2} & \cdots & a_{n-1} & a_{n}\\
a_{1} & x & a_{2} & \cdots & a_{n-1} & a_{n}\\
a_{1} & a_{2} & x & \cdots & a_{n-1} & a_{n}\\
\vdots & \vdots & \vdots & \ddots & \vdots & \vdots\\
a_{1} & a_{2} & a_{3} & \cdots & x & a_{n}\\
a_{1} & a_{2} & a_{3} & \cdots & a_{n} & x
\end{array}
\right)  =\left(  x+\sum_{i=1}^{n}a_{i}\right)  \prod_{i=1}^{n}\left(
x-a_{i}\right)  .
\]
This solves Exercise \ref{exe.det.a1a2anx}.
\end{proof}
\end{verlong}

\subsection{Solution to Exercise \ref{exe.det.2diags}}

\begin{vershort}
\begin{proof}
[Solution to Exercise \ref{exe.det.2diags}.]Let $z$ be the permutation
$\operatorname*{cyc}\nolimits_{n,n-1,\ldots,1}$ (where we are using the
notations of Definition \ref{def.perm.cycles}). Then, every $i\in\left\{
1,2,\ldots,n\right\}  $ satisfies%
\begin{equation}
z\left(  i\right)  =%
\begin{cases}
i-1, & \text{if }i>1;\\
n, & \text{if }i=1
\end{cases}
. \label{sol.det.2diags.short.z}%
\end{equation}
(This follows easily from the definition of $z$.) In particular, $z\left(
1\right)  =n\neq1$ (since $n>1$), and thus $z\neq\operatorname*{id}$. Every
$i\in\left\{  1,2,\ldots,n\right\}  $ satisfies%
\begin{align}
z\left(  i\right)   &  =%
\begin{cases}
i-1, & \text{if }i>1;\\
n, & \text{if }i=1
\end{cases}
\nonumber\\
&  \equiv%
\begin{cases}
i-1, & \text{if }i>1;\\
i-1, & \text{if }i=1
\end{cases}
\ \ \ \ \ \ \ \ \ \ \left(
\begin{array}
[c]{c}%
\text{since }n\equiv0=\underbrace{1}_{=i}-1=i-1\operatorname{mod}n\\
\text{in the case when }i=1
\end{array}
\right) \nonumber\\
&  =i-1\operatorname{mod}n. \label{sol.det.2diags.short.zmod}%
\end{align}
Notice also that $z=\operatorname*{cyc}\nolimits_{n,n-1,\ldots,1}$, so that
$\left(  -1\right)  ^{z}=\left(  -1\right)  ^{\operatorname*{cyc}%
\nolimits_{n,n-1,\ldots,1}}=\left(  -1\right)  ^{n-1}$ (by Exercise
\ref{exe.perm.cycles} \textbf{(d)}, applied to $k=n$ and $\left(  i_{1}%
,i_{2},\ldots,i_{k}\right)  =\left(  n,n-1,\ldots,1\right)  $).

We recall the following simple fact: If $p$ and $q$ are two elements of
$\left\{  1,2,\ldots,n\right\}  $ such that $p\equiv q\operatorname{mod}n$,
then $p=q$. We shall use this fact several times (tacitly) in the following arguments.

Now, I claim that if $\sigma\in S_{n}$ satisfies $\sigma\notin\left\{
\operatorname*{id},z\right\}  $, then%
\begin{equation}
\left(
\begin{array}
[c]{c}%
\text{there exists an }i\in\left\{  1,2,\ldots,n\right\}  \text{ satisfying}\\
i\neq\sigma\left(  i\right)  \text{ and }i\not \equiv \sigma\left(  i\right)
+1\operatorname{mod}n
\end{array}
\right)  . \label{sol.det.2diags.short.exi}%
\end{equation}

[\textit{Proof of (\ref{sol.det.2diags.short.exi}):} Let $\sigma\in S_{n}$ be
such that $\sigma\notin\left\{  \operatorname*{id},z\right\}  $. Thus,
$\sigma\neq\operatorname*{id}$ and $\sigma\neq z$.

We need to prove (\ref{sol.det.2diags.short.exi}). Indeed, let us assume the
contrary (for the sake of contradiction). Thus,%
\begin{equation}
\text{every }i\in\left\{  1,2,\ldots,n\right\}  \text{ satisfies either
}i=\sigma\left(  i\right)  \text{ or }i\equiv\sigma\left(  i\right)
+1\operatorname{mod}n. \label{sol.det.2diags.short.exi.pf.ass}%
\end{equation}
\footnote{I use the words \textquotedblleft either\textquotedblleft%
/\textquotedblleft or\textquotedblright\ in a non-exclusive meaning (i.e.,
when I say \textquotedblleft either $\mathcal{A}$ or $\mathcal{B}%
$\textquotedblright, I mean to include also the case when both $\mathcal{A}$
and $\mathcal{B}$ hold simultaneously), but here it does not matter (because
we cannot have $i=\sigma\left(  i\right)  $ and $i\equiv\sigma\left(
i\right)  +1\operatorname{mod}n$ at the same time).}

There exists a $J\in\left\{  1,2,\ldots,n\right\}  $ such that $\sigma\left(
J\right)  \neq J$ (since $\sigma\neq\operatorname*{id}$). Let $j$ be the
smallest such $J$. Thus, $\sigma\left(  j\right)  \neq j$, but every $J<j$
satisfies $\sigma\left(  J\right)  =J$.

Applying (\ref{sol.det.2diags.short.exi.pf.ass}) to $i=j$, we see that either
$j=\sigma\left(  j\right)  $ or $j\equiv\sigma\left(  j\right)
+1\operatorname{mod}n$. Since $j=\sigma\left(  j\right)  $ cannot hold
(because we have $\sigma\left(  j\right)  \neq j$), we thus have
$j\equiv\sigma\left(  j\right)  +1\operatorname{mod}n$. In other words,
$\sigma\left(  j\right)  \equiv j-1\operatorname{mod}n$.

We have $j=1$\ \ \ \ \footnote{\textit{Proof.} Assume the contrary. Thus,
$j\neq1$, so that $j>1$. Hence, $j-1\in\left\{  1,2,\ldots,n\right\}  $. Also,
$j-1<j$. Hence, $\sigma\left(  j-1\right)  =j-1$ (since every $J<j$ satisfies
$\sigma\left(  J\right)  =J$). Thus, $\sigma\left(  j-1\right)  =j-1\equiv
\sigma\left(  j\right)  \operatorname{mod}n$. Since both $\sigma\left(
j-1\right)  $ and $\sigma\left(  j\right)  $ are elements of $\left\{
1,2,\ldots,n\right\}  $, this shows that $\sigma\left(  j-1\right)
=\sigma\left(  j\right)  $, and thus $j-1=j$ (since $\sigma$ is injective).
But this is absurd. This contradiction shows that our assumption was wrong,
qed.}. Thus, $\sigma\left(  \underbrace{1}_{=j}\right)  =\sigma\left(
j\right)  \equiv\underbrace{j}_{=1}-1=0\equiv n\operatorname{mod}n$. Since
both $\sigma\left(  1\right)  $ and $n$ belong to $\left\{  1,2,\ldots
,n\right\}  $, this shows that $\sigma\left(  1\right)  =n$.

There exists a $K\in\left\{  1,2,\ldots,n\right\}  $ such that $\sigma\left(
K\right)  =K$\ \ \ \ \footnote{\textit{Proof.} We have $\sigma\neq z$. Hence,
there exists an $i\in\left\{  1,2,\ldots,n\right\}  $ such that $\sigma\left(
i\right)  \neq z\left(  i\right)  $. Consider this $i$.
\par
If we had $i\equiv\sigma\left(  i\right)  +1\operatorname{mod}n$, then we
would have $\sigma\left(  i\right)  \equiv i-1\equiv z\left(  i\right)
\operatorname{mod}n$ (by (\ref{sol.det.2diags.short.zmod})), which would
entail $\sigma\left(  i\right)  =z\left(  i\right)  $ (since both
$\sigma\left(  i\right)  $ and $z\left(  i\right)  $ belong to $\left\{
1,2,\ldots,n\right\}  $); but this would contradict $\sigma\left(  i\right)
\neq z\left(  i\right)  $. Hence, we cannot have $i\equiv\sigma\left(
i\right)  +1\operatorname{mod}n$.
\par
We have either $i=\sigma\left(  i\right)  $ or $i\equiv\sigma\left(  i\right)
+1\operatorname{mod}n$ (because of (\ref{sol.det.2diags.short.exi.pf.ass})).
Thus, $i=\sigma\left(  i\right)  $ (since we cannot have $i\equiv\sigma\left(
i\right)  +1\operatorname{mod}n$). Hence, there exists a $K\in\left\{
1,2,\ldots,n\right\}  $ such that $\sigma\left(  K\right)  =K$ (namely,
$K=i$). Qed.}. Let $k$ be the largest such $K$. Thus, $\sigma\left(  k\right)
=k$, but every $K>k$ satisfies $\sigma\left(  K\right)  \neq K$.

We have $n\neq1$ and thus $\sigma\left(  n\right)  \neq\sigma\left(  1\right)
$ (since $\sigma$ is injective). Thus, $\sigma\left(  n\right)  \neq
\sigma\left(  1\right)  =n$.

We cannot have $k=n$ (because otherwise, we would have $\sigma\left(
\underbrace{n}_{=k}\right)  =\sigma\left(  k\right)  =k=n$, which would
contradict $\sigma\left(  n\right)  \neq n$). Thus, $k<n$. Hence,
$k+1\in\left\{  1,2,\ldots,n\right\}  $. Therefore, $\sigma\left(  k+1\right)
\neq k+1$ (since every $K>k$ satisfies $\sigma\left(  K\right)  \neq K$, and
since $k+1>k$). Now, applying (\ref{sol.det.2diags.short.exi.pf.ass}) to
$i=k+1$, we conclude that either $k+1=\sigma\left(  k+1\right)  $ or
$k+1\equiv\sigma\left(  k+1\right)  +1\operatorname{mod}n$. Since we cannot
have $k+1=\sigma\left(  k+1\right)  $ (because $\sigma\left(  k+1\right)  \neq
k+1$), we thus must have $k+1\equiv\sigma\left(  k+1\right)
+1\operatorname{mod}n$. In other words, $k\equiv\sigma\left(  k+1\right)
\operatorname{mod}n$. Hence, $k=\sigma\left(  k+1\right)  $ (since both $k$
and $\sigma\left(  k+1\right)  $ belong to $\left\{  1,2,\ldots,n\right\}  $),
so that $\sigma\left(  k+1\right)  =k=\sigma\left(  k\right)  $. Since
$\sigma$ is injective, this yields $k+1=k$, which is absurd. This
contradiction shows that our assumption was wrong. Hence,
(\ref{sol.det.2diags.short.exi}) is proven.]

Let us now write our matrix $A$ in the form $A=\left(  a_{i,j}\right)  _{1\leq
i\leq n,\ 1\leq j\leq n}$. Then,%
\[
\left(  a_{i,j}\right)  _{1\leq i\leq n,\ 1\leq j\leq n}=A=\left(
\begin{cases}
a_{j}, & \text{if }i=j;\\
b_{j}, & \text{if }i\equiv j+1\operatorname{mod}n;\\
0, & \text{otherwise}%
\end{cases}
\right)  _{1\leq i\leq n,\ 1\leq j\leq n}.
\]
In other words, we have%
\begin{equation}
a_{i,j}=%
\begin{cases}
a_{j}, & \text{if }i=j;\\
b_{j}, & \text{if }i\equiv j+1\operatorname{mod}n;\\
0, & \text{otherwise}%
\end{cases}
\label{sol.det.2diags.short.aij}%
\end{equation}
for every $\left(  i,j\right)  \in\left\{  1,2,\ldots,n\right\}  ^{2}$.

For every $i\in\left\{  1,2,\ldots,n\right\}  $, we have%
\begin{align}
a_{i,i}  &  =%
\begin{cases}
a_{i}, & \text{if }i=i;\\
b_{i}, & \text{if }i\equiv i+1\operatorname{mod}n;\\
0, & \text{otherwise}%
\end{cases}
\ \ \ \ \ \ \ \ \ \ \left(  \text{by (\ref{sol.det.2diags.short.aij}), applied
to }\left(  i,i\right)  \text{ instead of }\left(  i,j\right)  \right)
\nonumber\\
&  =a_{i}\ \ \ \ \ \ \ \ \ \ \left(  \text{since }i=i\right)
\label{sol.det.2diags.short.aii}%
\end{align}
and%
\begin{align}
a_{i,z\left(  i\right)  }  &  =%
\begin{cases}
a_{z\left(  i\right)  }, & \text{if }i=z\left(  i\right)  ;\\
b_{z\left(  i\right)  }, & \text{if }i\equiv z\left(  i\right)
+1\operatorname{mod}n;\\
0, & \text{otherwise}%
\end{cases}
\nonumber\\
&  \ \ \ \ \ \ \ \ \ \ \left(  \text{by (\ref{sol.det.2diags.short.aij}),
applied to }\left(  i,z\left(  i\right)  \right)  \text{ instead of }\left(
i,j\right)  \right) \nonumber\\
&  =b_{z\left(  i\right)  }\ \ \ \ \ \ \ \ \ \ \left(  \text{since }i\equiv
z\left(  i\right)  +1\operatorname{mod}n\text{ (by
(\ref{sol.det.2diags.short.zmod}))}\right)  .
\label{sol.det.2diags.short.aizi}%
\end{align}

It is now easy to see that if $\sigma\in S_{n}$ satisfies $\sigma
\notin\left\{  \operatorname*{id},z\right\}  $, then
\begin{equation}
\prod_{i=1}^{n}a_{i,\sigma\left(  i\right)  }=0 \label{sol.det.2diags.short.0}%
\end{equation}
\footnote{\textit{Proof of (\ref{sol.det.2diags.short.0}):} Let $\sigma\in
S_{n}$ be such that $\sigma\notin\left\{  \operatorname*{id},z\right\}  $.
Thus, there exists an $i\in\left\{  1,2,\ldots,n\right\}  $ satisfying
$i\neq\sigma\left(  i\right)  $ and $i\not \equiv \sigma\left(  i\right)
+1\operatorname{mod}n$ (because of (\ref{sol.det.2diags.short.exi})). Consider
this $i$.
\par
From (\ref{sol.det.2diags.short.aij}) (applied to $\left(  i,\sigma\left(
i\right)  \right)  $ instead of $\left(  i,j\right)  $), we obtain%
\[
a_{i,\sigma\left(  i\right)  }=%
\begin{cases}
a_{\sigma\left(  i\right)  }, & \text{if }i=\sigma\left(  i\right)  ;\\
b_{\sigma\left(  i\right)  }, & \text{if }i\equiv\sigma\left(  i\right)
+1\operatorname{mod}n;\\
0, & \text{otherwise}%
\end{cases}
=0
\]
(since $i\neq\sigma\left(  i\right)  $ and $i\not \equiv \sigma\left(
i\right)  +1\operatorname{mod}n$).
\par
Now, let us forget that we fixed $i$. We thus have shown that there exists an
$i\in\left\{  1,2,\ldots,n\right\}  $ satisfying $a_{i,\sigma\left(  i\right)
}=0$. In other words, at least one factor of the product $\prod_{i=1}%
^{n}a_{i,\sigma\left(  i\right)  }$ equals $0$. Therefore, the whole product
$\prod_{i=1}^{n}a_{i,\sigma\left(  i\right)  }$ equals $0$. This proves
(\ref{sol.det.2diags.short.0}).}.

Now, (\ref{eq.det.eq.2}) becomes%
\begin{align*}
\det A  &  =\sum_{\sigma\in S_{n}}\left(  -1\right)  ^{\sigma}\prod_{i=1}%
^{n}a_{i,\sigma\left(  i\right)  }\\
&  =\underbrace{\left(  -1\right)  ^{\operatorname*{id}}}_{=1}\prod_{i=1}%
^{n}\underbrace{a_{i,\operatorname*{id}\left(  i\right)  }}_{=a_{i,i}%
}+\underbrace{\left(  -1\right)  ^{z}}_{=\left(  -1\right)  ^{n-1}}\prod
_{i=1}^{n}a_{i,z\left(  i\right)  }+\sum_{\substack{\sigma\in S_{n}%
;\\\sigma\notin\left\{  \operatorname*{id},z\right\}  }}\left(  -1\right)
^{\sigma}\underbrace{\prod_{i=1}^{n}a_{i,\sigma\left(  i\right)  }%
}_{\substack{=0\\\text{(by (\ref{sol.det.2diags.short.0}))}}}\\
&  \ \ \ \ \ \ \ \ \ \ \left(
\begin{array}
[c]{c}%
\text{here, we have split off the addends for }\sigma=\operatorname*{id}\text{
and}\\
\text{for }\sigma=z\text{ from the sum (since }z\neq\operatorname*{id}\text{)}%
\end{array}
\right) \\
&  =\prod_{i=1}^{n}a_{i,i}+\left(  -1\right)  ^{n-1}\prod_{i=1}^{n}%
a_{i,z\left(  i\right)  }+\underbrace{\sum_{\substack{\sigma\in S_{n}%
;\\\sigma\notin\left\{  \operatorname*{id},z\right\}  }}\left(  -1\right)
^{\sigma}0}_{=0}\\
&  =\prod_{i=1}^{n}\underbrace{a_{i,i}}_{\substack{=a_{i}\\\text{(by
(\ref{sol.det.2diags.short.aii}))}}}+\left(  -1\right)  ^{n-1}\prod_{i=1}%
^{n}\underbrace{a_{i,z\left(  i\right)  }}_{\substack{=b_{z\left(  i\right)
}\\\text{(by (\ref{sol.det.2diags.short.aizi}))}}}\\
&  =\prod_{i=1}^{n}a_{i}+\left(  -1\right)  ^{n-1}\underbrace{\prod_{i=1}%
^{n}b_{z\left(  i\right)  }}_{\substack{=\prod_{i=1}^{n}b_{i}\\\text{(here, we
have substituted }i\\\text{for }z\left(  i\right)  \text{ in the
product,}\\\text{since }z:\left\{  1,2,\ldots,n\right\}  \rightarrow\left\{
1,2,\ldots,n\right\}  \\\text{is a bijection)}}}=\underbrace{\prod_{i=1}%
^{n}a_{i}}_{=a_{1}a_{2}\cdots a_{n}}+\left(  -1\right)  ^{n-1}%
\underbrace{\prod_{i=1}^{n}b_{i}}_{=b_{1}b_{2}\cdots b_{n}}\\
&  =a_{1}a_{2}\cdots a_{n}+\left(  -1\right)  ^{n-1}b_{1}b_{2}\cdots b_{n}.
\end{align*}
This solves Exercise \ref{exe.det.2diags}.
\end{proof}
\end{vershort}

\begin{verlong}
\begin{proof}
[Solution to Exercise \ref{exe.det.2diags}.]We first make some preliminary definitions.

Let $\left[  n\right]  $ denote the set $\left\{  1,2,\ldots,n\right\}  $. The
set $S_{n}$ is the set of all permutations of $\left\{  1,2,\ldots,n\right\}
$. In other words, the set $S_{n}$ is the set of all permutations of $\left[
n\right]  $ (since $\left[  n\right]  =\left\{  1,2,\ldots,n\right\}  $).

The integers $1,2,\ldots,n$ take all possible remainders when divided by $n$,
and take each of them exactly once. In other words: For every $h\in\mathbb{Z}%
$, there exists a unique element $g\in\left\{  1,2,\ldots,n\right\}  $
satisfying $g\equiv h\operatorname{mod}n$\ \ \ \ \footnote{\textit{Proof.} Let
$h\in\mathbb{Z}$. Exercise \ref{exe.mod.unique-cong} (applied to $N=n$ and
$p=0$) shows that there exists a unique element $g\in\left\{  0+1,0+2,\ldots
,0+n\right\}  $ satisfying $g\equiv h\operatorname{mod}n$. In other words,
there exists a unique element $g\in\left\{  1,2,\ldots,n\right\}  $ satisfying
$g\equiv h\operatorname{mod}n$ (since $\left\{  0+1,0+2,\ldots,0+n\right\}
=\left\{  1,2,\ldots,n\right\}  $). Qed.}. Since $\left\{  1,2,\ldots
,n\right\}  =\left[  n\right]  $, this rewrites as follows: For every
$h\in\mathbb{Z}$, there exists a unique element $g\in\left[  n\right]  $
satisfying $g\equiv h\operatorname{mod}n$. We shall denote this $g$ by
$\operatorname*{posrem}h$. Thus, we have the following facts:

\begin{itemize}
\item For every $h\in\mathbb{Z}$, we have%
\begin{equation}
\operatorname*{posrem}h\in\left[  n\right]  \ \ \ \ \ \ \ \ \ \ \text{and}%
\ \ \ \ \ \ \ \ \ \ \operatorname*{posrem}h\equiv h\operatorname{mod}n
\label{sol.det.2diags.posrem.1}%
\end{equation}
\footnote{\textit{Proof of (\ref{sol.det.2diags.posrem.1}):} Let
$h\in\mathbb{Z}$. We know that $\operatorname*{posrem}h$ is the unique element
$g\in\left[  n\right]  $ satisfying $g\equiv h\operatorname{mod}n$ (because
this is how we have defined $\operatorname*{posrem}h$). Thus,
$\operatorname*{posrem}h$ is an element $g\in\left[  n\right]  $ satisfying
$g\equiv h\operatorname{mod}n$. In other words, $\operatorname*{posrem}h$ is
an element of $\left[  n\right]  $ and satisfies $\operatorname*{posrem}%
h\equiv h\operatorname{mod}n$. In other words, $\operatorname*{posrem}%
h\in\left[  n\right]  $ and $\operatorname*{posrem}h\equiv h\operatorname{mod}%
n$. This proves (\ref{sol.det.2diags.posrem.1}).}.

\item If $h\in\mathbb{Z}$ and $k\in\left[  n\right]  $ are such that $k\equiv
h\operatorname{mod}n$, then%
\begin{equation}
k=\operatorname*{posrem}h \label{sol.det.2diags.posrem.2}%
\end{equation}
\footnote{\textit{Proof of (\ref{sol.det.2diags.posrem.2}):} Let
$h\in\mathbb{Z}$ and $k\in\left[  n\right]  $ be such that $k\equiv
h\operatorname{mod}n$.
\par
We know that $\operatorname*{posrem}h$ is the unique element $g\in\left[
n\right]  $ satisfying $g\equiv h\operatorname{mod}n$ (because this is how we
have defined $\operatorname*{posrem}h$). Thus, in particular,
$\operatorname*{posrem}h$ is the only such element. In other words, if
$g\in\left[  n\right]  $ is any element satisfying $g\equiv
h\operatorname{mod}n$, then $g=\operatorname*{posrem}h$. Applying this to
$g=k$, we obtain $k=\operatorname*{posrem}h$ (since $k\in\left[  n\right]  $
is an element satisfying $k\equiv h\operatorname{mod}n$). This proves
(\ref{sol.det.2diags.posrem.2}).}.

\item If $a\in\left[  n\right]  $ and $b\in\left[  n\right]  $ are such that
$a\equiv b\operatorname{mod}n$, then%
\begin{equation}
a=b \label{sol.det.2diags.posrem.3}%
\end{equation}
\footnote{\textit{Proof of (\ref{sol.det.2diags.posrem.3}):} Let $a\in\left[
n\right]  $ and $b\in\left[  n\right]  $ be such that $a\equiv
b\operatorname{mod}n$. Applying (\ref{sol.det.2diags.posrem.2}) to $h=b$ and
$k=a$, we obtain $a=\operatorname*{posrem}b$ (since $a\equiv
b\operatorname{mod}n$). Applying (\ref{sol.det.2diags.posrem.2}) to $h=b$ and
$k=b$, we obtain $b=\operatorname*{posrem}b$ (since $b\equiv
b\operatorname{mod}n$). Thus, $a=\operatorname*{posrem}b=b$. This proves
(\ref{sol.det.2diags.posrem.3}).}.
\end{itemize}

Now, fix $a\in\mathbb{Z}$. Every $g\in\left[  n\right]  $ satisfies
$\operatorname*{posrem}\left(  g+a\right)  \in\left[  n\right]  $%
\ \ \ \ \footnote{\textit{Proof.} Let $g\in\left[  n\right]  $. Then,
(\ref{sol.det.2diags.posrem.1}) (applied to $h=g+a$) shows that
$\operatorname*{posrem}\left(  g+a\right)  \in\left[  n\right]  $ and
$\operatorname*{posrem}\left(  g+a\right)  \equiv g+a\operatorname{mod}n$. In
particular, $\operatorname*{posrem}\left(  g+a\right)  \in\left[  n\right]  $,
qed.}. Hence, we can define a map $\operatorname*{shift}\nolimits_{a}:\left[
n\right]  \rightarrow\left[  n\right]  $ by setting%
\[
\left(  \operatorname*{shift}\nolimits_{a}\left(  g\right)
=\operatorname*{posrem}\left(  g+a\right)  \ \ \ \ \ \ \ \ \ \ \text{for every
}g\in\left[  n\right]  \right)  .
\]
Consider this map $\operatorname*{shift}\nolimits_{a}$.

Let us now forget that we fixed $a\in\mathbb{Z}$. We thus have defined a map
$\operatorname*{shift}\nolimits_{a}:\left[  n\right]  \rightarrow\left[
n\right]  $ for every $a\in\mathbb{Z}$. (For example, if $n=5$ and $a=2$, then
the map $\operatorname*{shift}\nolimits_{a}$ sends $1,2,3,4,5$ to $3,4,5,1,2$, respectively.)

We observe the following facts:

\begin{itemize}
\item We have%
\begin{equation}
\operatorname*{shift}\nolimits_{0}=\operatorname*{id}
\label{sol.det.2diags.shift.0}%
\end{equation}
\footnote{\textit{Proof of (\ref{sol.det.2diags.shift.0}):} Let $g\in\left[
n\right]  $. Then, $g\equiv g\operatorname{mod}n$. Hence,
(\ref{sol.det.2diags.posrem.2}) (applied to $h=g$ and $k=g$) shows that
$g=\operatorname*{posrem}g$. Now, the definition of $\operatorname*{shift}%
\nolimits_{0}$ shows that $\operatorname*{shift}\nolimits_{0}\left(  g\right)
=\operatorname*{posrem}\left(  \underbrace{g+0}_{=g}\right)
=\operatorname*{posrem}g$. Compared with $\operatorname*{id}\left(  g\right)
=g=\operatorname*{posrem}g$, this yields $\operatorname*{shift}\nolimits_{0}%
\left(  g\right)  =\operatorname*{id}\left(  g\right)  $.
\par
Now, let us forget that we fixed $g$. We thus have proven that
$\operatorname*{shift}\nolimits_{0}\left(  g\right)  =\operatorname*{id}%
\left(  g\right)  $ for every $g\in\left[  n\right]  $. In other words,
$\operatorname*{shift}\nolimits_{0}=\operatorname*{id}$. This proves
(\ref{sol.det.2diags.shift.0}).}.

\item We have%
\begin{equation}
\operatorname*{shift}\nolimits_{a+b}=\operatorname*{shift}\nolimits_{a}%
\circ\operatorname*{shift}\nolimits_{b}\ \ \ \ \ \ \ \ \ \ \text{for any }%
a\in\mathbb{Z}\text{ and }b\in\mathbb{Z} \label{sol.det.2diags.shift.a+b}%
\end{equation}
\footnote{\textit{Proof of (\ref{sol.det.2diags.shift.a+b}):} Let
$a\in\mathbb{Z}$ and $b\in\mathbb{Z}$.
\par
Let $g\in\left[  n\right]  $. Then, $\operatorname*{shift}\nolimits_{a+b}%
\left(  g\right)  =\operatorname*{posrem}\left(  g+a+b\right)  $ (by the
definition of $\operatorname*{shift}\nolimits_{a+b}$).
\par
On the other hand, $\operatorname*{shift}\nolimits_{b}\left(  g\right)
=\operatorname*{posrem}\left(  g+b\right)  $ (by the definition of
$\operatorname*{shift}\nolimits_{b}$). Applying (\ref{sol.det.2diags.posrem.1}%
) to $h=g+b$, we obtain $\operatorname*{posrem}\left(  g+b\right)  \in\left[
n\right]  $ and $\operatorname*{posrem}\left(  g+b\right)  \equiv
g+b\operatorname{mod}n$.
\par
Moreover, $\left(  \operatorname*{shift}\nolimits_{a}\circ
\operatorname*{shift}\nolimits_{b}\right)  \left(  g\right)
=\operatorname*{shift}\nolimits_{a}\left(  \operatorname*{shift}%
\nolimits_{b}\left(  g\right)  \right)  =\operatorname*{posrem}\left(
\operatorname*{shift}\nolimits_{b}\left(  g\right)  +a\right)  $ (by the
definition of $\operatorname*{shift}\nolimits_{a}$). Applying
(\ref{sol.det.2diags.posrem.1}) to $h=\operatorname*{shift}\nolimits_{b}%
\left(  g\right)  +a$, we obtain $\operatorname*{posrem}\left(
\operatorname*{shift}\nolimits_{b}\left(  g\right)  +a\right)  \in\left[
n\right]  $ and $\operatorname*{posrem}\left(  \operatorname*{shift}%
\nolimits_{b}\left(  g\right)  +a\right)  \equiv\operatorname*{shift}%
\nolimits_{b}\left(  g\right)  +a\operatorname{mod}n$.
\par
Now,%
\begin{align*}
&  \left(  \operatorname*{shift}\nolimits_{a}\circ\operatorname*{shift}%
\nolimits_{b}\right)  \left(  g\right) \\
&  =\operatorname*{posrem}\left(  \operatorname*{shift}\nolimits_{b}\left(
g\right)  +a\right)  \equiv\underbrace{\operatorname*{shift}\nolimits_{b}%
\left(  g\right)  }_{=\operatorname*{posrem}\left(  g+b\right)  \equiv
g+b\operatorname{mod}n}+a\\
&  \equiv g+b+a=g+a+b\operatorname{mod}n.
\end{align*}
Thus, we can apply (\ref{sol.det.2diags.posrem.3}) to $h=g+a+b$ and $k=\left(
\operatorname*{shift}\nolimits_{a}\circ\operatorname*{shift}\nolimits_{b}%
\right)  \left(  g\right)  $ (since $\left(  \operatorname*{shift}%
\nolimits_{a}\circ\operatorname*{shift}\nolimits_{b}\right)  \left(  g\right)
=\operatorname*{posrem}\left(  \operatorname*{shift}\nolimits_{b}\left(
g\right)  +a\right)  \in\left[  n\right]  $). As a result, we obtain $\left(
\operatorname*{shift}\nolimits_{a}\circ\operatorname*{shift}\nolimits_{b}%
\right)  \left(  g\right)  =\operatorname*{posrem}\left(  g+a+b\right)  $.
Comparing this with $\operatorname*{shift}\nolimits_{a+b}\left(  g\right)
=\operatorname*{posrem}\left(  g+a+b\right)  $, we obtain $\left(
\operatorname*{shift}\nolimits_{a}\circ\operatorname*{shift}\nolimits_{b}%
\right)  \left(  g\right)  =\operatorname*{shift}\nolimits_{a+b}\left(
g\right)  $.
\par
Now, let us forget that we fixed $g$. We thus have shown that $\left(
\operatorname*{shift}\nolimits_{a}\circ\operatorname*{shift}\nolimits_{b}%
\right)  \left(  g\right)  =\operatorname*{shift}\nolimits_{a+b}\left(
g\right)  $ for every $g\in\left[  n\right]  $. In other words,
$\operatorname*{shift}\nolimits_{a}\circ\operatorname*{shift}\nolimits_{b}%
=\operatorname*{shift}\nolimits_{a+b}$. This proves
(\ref{sol.det.2diags.shift.a+b}).}.

\item We have%
\begin{equation}
\operatorname*{shift}\nolimits_{a}\in S_{n}\ \ \ \ \ \ \ \ \ \ \text{for every
}a\in\mathbb{Z} \label{sol.det.2diags.shift.Sn}%
\end{equation}
\footnote{\textit{Proof of (\ref{sol.det.2diags.shift.Sn}):} Let
$a\in\mathbb{Z}$. From (\ref{sol.det.2diags.shift.a+b}) (applied to $b=-a$),
we obtain $\operatorname*{shift}\nolimits_{a+\left(  -a\right)  }%
=\operatorname*{shift}\nolimits_{a}\circ\operatorname*{shift}\nolimits_{-a}$.
Hence,%
\begin{align*}
\operatorname*{shift}\nolimits_{a}\circ\operatorname*{shift}\nolimits_{-a}  &
=\operatorname*{shift}\nolimits_{a+\left(  -a\right)  }=\operatorname*{shift}%
\nolimits_{0}\ \ \ \ \ \ \ \ \ \ \left(  \text{since }a+\left(  -a\right)
=0\right) \\
&  =\operatorname*{id}\ \ \ \ \ \ \ \ \ \ \left(  \text{by
(\ref{sol.det.2diags.shift.0})}\right)  .
\end{align*}
\par
From (\ref{sol.det.2diags.shift.a+b}) (applied to $-a$ and $a$ instead of $a$
and $b$), we obtain $\operatorname*{shift}\nolimits_{\left(  -a\right)
+a}=\operatorname*{shift}\nolimits_{-a}\circ\operatorname*{shift}%
\nolimits_{a}$. Hence,%
\begin{align*}
\operatorname*{shift}\nolimits_{-a}\circ\operatorname*{shift}\nolimits_{a}  &
=\operatorname*{shift}\nolimits_{\left(  -a\right)  +a}=\operatorname*{shift}%
\nolimits_{0}\ \ \ \ \ \ \ \ \ \ \left(  \text{since }\left(  -a\right)
+a=0\right) \\
&  =\operatorname*{id}\ \ \ \ \ \ \ \ \ \ \left(  \text{by
(\ref{sol.det.2diags.shift.0})}\right)  .
\end{align*}
\par
The maps $\operatorname*{shift}\nolimits_{a}$ and $\operatorname*{shift}%
\nolimits_{-a}$ are mutually inverse (since $\operatorname*{shift}%
\nolimits_{a}\circ\operatorname*{shift}\nolimits_{-a}=\operatorname*{id}$ and
$\operatorname*{shift}\nolimits_{-a}\circ\operatorname*{shift}\nolimits_{a}%
=\operatorname*{id}$). Thus, the map $\operatorname*{shift}\nolimits_{a}$ is
invertible, and therefore a bijection.
\par
So the map $\operatorname*{shift}\nolimits_{a}$ is a bijection $\left[
n\right]  \rightarrow\left[  n\right]  $. In other words, the map
$\operatorname*{shift}\nolimits_{a}$ is a permutation of the set $\left[
n\right]  $. In other words, the map $\operatorname*{shift}\nolimits_{a}$ is a
permutation of the set $\left\{  1,2,\ldots,n\right\}  $ (since $\left[
n\right]  =\left\{  1,2,\ldots,n\right\}  $). In other words,
$\operatorname*{shift}\nolimits_{a}\in S_{n}$ (since $S_{n}$ is the set of all
permutations of the set $\left\{  1,2,\ldots,n\right\}  $). This proves
(\ref{sol.det.2diags.shift.Sn}).}.
\end{itemize}

Let $z=\operatorname*{shift}\nolimits_{-1}$. Then, $z=\operatorname*{shift}%
\nolimits_{-1}\in S_{n}$ (by (\ref{sol.det.2diags.shift.Sn}), applied to
$a=-1$). In other words, $z$ is a permutation of $\left[  n\right]  $ (since
$S_{n}$ is the set of all permutations of $\left[  n\right]  $). We observe
some properties of $z$:

\begin{itemize}
\item We have $z=\operatorname*{cyc}\nolimits_{n,n-1,\ldots,1}$, where we are
using the notations of Definition \ref{def.perm.cycles}. We will not use this
fact, and thus we will not prove it (but the proof is almost trivial).

\item We have%
\begin{equation}
z\left(  i\right)  \equiv i-1\operatorname{mod}n\ \ \ \ \ \ \ \ \ \ \text{for
every }i\in\left[  n\right]  \label{sol.det.2diags.zmod}%
\end{equation}
\footnote{\textit{Proof of (\ref{sol.det.2diags.zmod}):} Let $i\in\left[
n\right]  $. Applying (\ref{sol.det.2diags.posrem.1}) to $h=i-1$, we obtain
$\operatorname*{posrem}\left(  i-1\right)  \in\left[  n\right]  $ and
$\operatorname*{posrem}\left(  i-1\right)  \equiv i-1\operatorname{mod}n$.
Now,%
\begin{align*}
\underbrace{z}_{=\operatorname*{shift}\nolimits_{-1}}\left(  i\right)   &
=\operatorname*{shift}\nolimits_{-1}\left(  i\right)  =\operatorname*{posrem}%
\left(  \underbrace{i+\left(  -1\right)  }_{=i-1}\right)
\ \ \ \ \ \ \ \ \ \ \left(  \text{by the definition of }\operatorname*{shift}%
\nolimits_{-1}\right) \\
&  =\operatorname*{posrem}\left(  i-1\right)  \equiv i-1\operatorname{mod}n.
\end{align*}
This proves (\ref{sol.det.2diags.zmod}).}. In other words,%
\begin{equation}
z\left(  i\right)  +1\equiv i\operatorname{mod}n\ \ \ \ \ \ \ \ \ \ \text{for
every }i\in\left[  n\right]  . \label{sol.det.2diags.zmod2}%
\end{equation}
In particular,%
\begin{equation}
i\neq z\left(  i\right)  \ \ \ \ \ \ \ \ \ \ \text{for every }i\in\left[
n\right]  \label{sol.det.2diags.zmod0}%
\end{equation}
\footnote{\textit{Proof of (\ref{sol.det.2diags.zmod0}):} Let $i\in\left[
n\right]  $. We have $i\equiv z\left(  i\right)  +1\operatorname{mod}n$ (by
(\ref{sol.det.2diags.zmod2})), so that $i-z\left(  i\right)  \equiv
1\not \equiv 0\operatorname{mod}n$ (since $n>1$). In other words,
$i\not \equiv z\left(  i\right)  \operatorname{mod}n$. Hence, $i\neq z\left(
i\right)  $. This proves (\ref{sol.det.2diags.zmod0}).}.

\item We have
\begin{equation}
z\left(  1\right)  =n \label{sol.det.2diags.shift.z.1}%
\end{equation}
\footnote{\textit{Proof of (\ref{sol.det.2diags.shift.z.1}):} We have
$n\in\left\{  1,2,\ldots,n\right\}  =\left[  n\right]  $ and $n\equiv
0\operatorname{mod}n$. Hence, $n=\operatorname*{posrem}0$ (by
(\ref{sol.det.2diags.posrem.2}), applied to $k=n$ and $h=0$). Now,
\begin{align*}
\underbrace{z}_{=\operatorname*{shift}\nolimits_{-1}}\left(  1\right)   &
=\operatorname*{shift}\nolimits_{-1}\left(  1\right)  =\operatorname*{posrem}%
\underbrace{\left(  1+\left(  -1\right)  \right)  }_{=0}%
\ \ \ \ \ \ \ \ \ \ \left(  \text{by the definition of }\operatorname*{shift}%
\nolimits_{-1}\right) \\
&  =\operatorname*{posrem}0=n\ \ \ \ \ \ \ \ \ \ \left(  \text{since
}n=\operatorname*{posrem}0\right)  .
\end{align*}
This proves (\ref{sol.det.2diags.shift.z.1}).}.

\item We have%
\begin{equation}
z\left(  i\right)  =i-1 \label{sol.det.2diags.shift.z.i}%
\end{equation}
for every $i\in\left[  n\right]  $ satisfying $i\neq1$%
\ \ \ \ \footnote{\textit{Proof of (\ref{sol.det.2diags.shift.z.i}):} Let
$i\in\left[  n\right]  $ be such that $i\neq1$. We have $i\in\left[  n\right]
=\left\{  1,2,\ldots,n\right\}  $. Combining this with $i\neq1$, we obtain
$i\in\left\{  1,2,\ldots,n\right\}  \setminus\left\{  1\right\}  =\left\{
2,3,\ldots,n\right\}  $, so that $i-1\in\left\{  1,2,\ldots,n-1\right\}
\subseteq\left\{  1,2,\ldots,n\right\}  =\left[  n\right]  $.
\par
We have $i-1\in\left[  n\right]  $ and $i-1\equiv i-1\operatorname{mod}n$.
Hence, $i-1=\operatorname*{posrem}\left(  i-1\right)  $ (by
(\ref{sol.det.2diags.posrem.2}), applied to $k=i-1$ and $h=i-1$). Now,
\begin{align*}
\underbrace{z}_{=\operatorname*{shift}\nolimits_{-1}}\left(  i\right)   &
=\operatorname*{shift}\nolimits_{-1}\left(  i\right)  =\operatorname*{posrem}%
\underbrace{\left(  i+\left(  -1\right)  \right)  }_{=i-1}%
\ \ \ \ \ \ \ \ \ \ \left(  \text{by the definition of }\operatorname*{shift}%
\nolimits_{-1}\right) \\
&  =\operatorname*{posrem}\left(  i-1\right)  =i-1\ \ \ \ \ \ \ \ \ \ \left(
\text{since }i-1=\operatorname*{posrem}\left(  i-1\right)  \right)  .
\end{align*}
This proves (\ref{sol.det.2diags.shift.z.i}).}.
\end{itemize}

For every $\sigma\in S_{n}$, we let $\operatorname*{Inv}\left(  \sigma\right)
$ be the set of all inversions of the permutation $\sigma$. Thus, for every
$\sigma\in S_{n}$, we have%
\begin{align}
\ell\left(  \sigma\right)   &  =\left(  \text{the number of inversions of
}\sigma\right)  \ \ \ \ \ \ \ \ \ \ \left(  \text{by the definition of }%
\ell\left(  \sigma\right)  \right) \nonumber\\
&  =\left\vert \underbrace{\left(  \text{the set of all inversions of }%
\sigma\right)  }_{\substack{=\operatorname*{Inv}\left(  \sigma\right)
\\\text{(since }\operatorname*{Inv}\left(  \sigma\right)  \text{ is the set of
all inversions of }\sigma\text{)}}}\right\vert \nonumber\\
&  =\left\vert \operatorname*{Inv}\left(  \sigma\right)  \right\vert .
\label{sol.det.2diags.invl}%
\end{align}

Now, we continue stating properties of $z$:

\begin{itemize}
\item The set $\operatorname*{Inv}\left(  z\right)  $ is the set of all
inversions of $z$ (by the definition of $\operatorname*{Inv}\left(  z\right)
$). For every $u\in\left\{  2,3,\ldots,n\right\}  $, we have $\left(
1,u\right)  \in\operatorname*{Inv}\left(  z\right)  $%
\ \ \ \ \footnote{\textit{Proof.} Let $u\in\left\{  2,3,\ldots,n\right\}  $.
Thus, $2\leq u\leq n$. Now, $1<2\leq u$. Thus, $1\leq1<u\leq n$. Moreover,
$u\in\left\{  2,3,\ldots,n\right\}  \subseteq\left\{  1,2,\ldots,n\right\}
=\left[  n\right]  $ and $u\neq1$ (since $1<u$). Hence,
(\ref{sol.det.2diags.shift.z.i}) (applied to $i=u$) shows that $z\left(
u\right)  =u-1<u\leq n$. But (\ref{sol.det.2diags.shift.z.1}) shows that
$z\left(  1\right)  =n>z\left(  u\right)  $ (since $z\left(  u\right)  <n$).
\par
Thus, $\left(  1,u\right)  $ is a pair $\left(  i,j\right)  $ of integers
satisfying $1\leq i<j\leq n$ and $z\left(  i\right)  >z\left(  j\right)  $
(since $1\leq1<u\leq n$ and $z\left(  1\right)  >z\left(  u\right)  $). In
other words, $\left(  1,u\right)  $ is an inversion of $z$ (by the definition
of an \textquotedblleft inversion of $z$\textquotedblright). In other words,
$\left(  1,u\right)  \in\operatorname*{Inv}\left(  z\right)  $ (since
$\operatorname*{Inv}\left(  z\right)  $ is the set of all inversions of $z$).
Qed.}. Hence, we can define a map $\rho:\left\{  2,3,\ldots,n\right\}
\rightarrow\operatorname*{Inv}\left(  z\right)  $ by%
\[
\left(  \rho\left(  u\right)  =\left(  1,u\right)
\ \ \ \ \ \ \ \ \ \ \text{for every }u\in\left\{  2,3,\ldots,n\right\}
\right)  .
\]
Consider this map $\rho$. The map $\rho$ is injective\footnote{\textit{Proof.}
Let $u_{1}$ and $u_{2}$ be two elements of $\left\{  2,3,\ldots,n\right\}  $
such that $\rho\left(  u_{1}\right)  =\rho\left(  u_{2}\right)  $. Then,
$\rho\left(  u_{1}\right)  =\left(  1,u_{1}\right)  $ (by the definition of
$\rho$) and $\rho\left(  u_{2}\right)  =\left(  1,u_{2}\right)  $ (by the
definition of $\rho$). Now, $\left(  1,u_{1}\right)  =\rho\left(
u_{1}\right)  =\rho\left(  u_{2}\right)  =\left(  1,u_{2}\right)  $. In other
words, $1=1$ and $u_{1}=u_{2}$.
\par
Let us now forget that we fixed $u_{1}$ and $u_{2}$. We thus have shown that
if $u_{1}$ and $u_{2}$ are two elements of $\left\{  2,3,\ldots,n\right\}  $
such that $\rho\left(  u_{1}\right)  =\rho\left(  u_{2}\right)  $, then
$u_{1}=u_{2}$. In other words, the map $\rho$ is injective, qed.} and
surjective\footnote{\textit{Proof.} Let $c\in\operatorname*{Inv}\left(
z\right)  $. Thus, $c$ is an element of the set $\operatorname*{Inv}\left(
z\right)  $. In other words, $c$ is an inversion of $z$ (since
$\operatorname*{Inv}\left(  z\right)  $ is the set of all inversions of $z$).
In other words, $c$ is a pair $\left(  i,j\right)  $ of integers satisfying
$1\leq i<j\leq n$ and $z\left(  i\right)  >z\left(  j\right)  $ (by the
definition of an \textquotedblleft inversion of $z$\textquotedblright).
Consider this pair $\left(  i,j\right)  $. Thus, $c=\left(  i,j\right)  $.
\par
We have $i\in\left[  n\right]  $ (since $1\leq i\leq n$) and $j\in\left[
n\right]  $ (since $1\leq j\leq n$). Also, $1\leq i<j$; thus, $j>1$, so that
$j\neq1$. Hence, $z\left(  j\right)  =j-1$ (by (\ref{sol.det.2diags.shift.z.i}%
), applied to $j$ instead of $i$). From $j>1$, we obtain $j\geq2$ (since $j$
is an integer). Hence, $j\in\left\{  2,3,\ldots,n\right\}  $ (since $j\leq
n$). Thus, $\rho\left(  j\right)  $ is well-defined. The definition of $\rho$
shows that $\rho\left(  j\right)  =\left(  1,j\right)  $.
\par
We assume (for the sake of contradiction) that $i\neq1$. Thus, $z\left(
i\right)  =i-1$ (by (\ref{sol.det.2diags.shift.z.i})). Hence, $z\left(
i\right)  =\underbrace{i}_{<j}-1<j-1=z\left(  j\right)  $, which contradicts
$z\left(  i\right)  >z\left(  j\right)  $. This contradiction shows that our
assumption (that $i\neq1$) was wrong. Hence, we cannot have $i\neq1$. We thus
have $i=1$. Now, $c=\left(  \underbrace{i}_{=1},j\right)  =\left(  1,j\right)
=\rho\left(  \underbrace{j}_{\in\left\{  2,3,\ldots,n\right\}  }\right)
\in\rho\left(  \left\{  2,3,\ldots,n\right\}  \right)  $.
\par
Let us now forget that we fixed $c$. We thus have shown that $c\in\rho\left(
\left\{  2,3,\ldots,n\right\}  \right)  $ for every $c\in\operatorname*{Inv}%
\left(  z\right)  $. In other words, $\operatorname*{Inv}\left(  z\right)
\subseteq\rho\left(  \left\{  2,3,\ldots,n\right\}  \right)  $. In other
words, the map $\rho$ is surjective, qed.}. Hence, the map $\rho$ is
bijective. Thus, we have found a bijective map between the sets $\left\{
2,3,\ldots,n\right\}  $ and $\operatorname*{Inv}\left(  z\right)  $ (namely,
$\rho$). Consequently,%
\[
\left\vert \operatorname*{Inv}\left(  z\right)  \right\vert =\left\vert
\left\{  2,3,\ldots,n\right\}  \right\vert =n-1.
\]
Now, (\ref{sol.det.2diags.invl}) (applied to $\sigma=z$) shows that
$\ell\left(  z\right)  =\left\vert \operatorname*{Inv}\left(  z\right)
\right\vert =n-1$. The definition of $\left(  -1\right)  ^{z}$ now shows that
\begin{equation}
\left(  -1\right)  ^{z}=\left(  -1\right)  ^{\ell\left(  z\right)  }=\left(
-1\right)  ^{n-1}\ \ \ \ \ \ \ \ \ \ \left(  \text{since }\ell\left(
z\right)  =n-1\right)  . \label{sol.det.2diags.signz}%
\end{equation}
(Alternatively, we could have obtained this equality from Exercise
\ref{exe.perm.cycles} \textbf{(d)} (applied to $k=n$ and $\left(  i_{1}%
,i_{2},\ldots,i_{k}\right)  =\left(  n,n-1,\ldots,1\right)  $) using the
observation that $z=\operatorname*{cyc}\nolimits_{n,n-1,\ldots,1}$.)

\item Now, let us recall that $n>1$; hence, $n\neq1$. We have $z\neq
\operatorname*{id}$\ \ \ \ \footnote{\textit{Proof.} Assume the contrary.
Thus, $z=\operatorname*{id}$. Hence, $\underbrace{z}_{=\operatorname*{id}%
}\left(  1\right)  =\operatorname*{id}\left(  1\right)  =1$. This contradicts
$z\left(  1\right)  =n\neq1$. This contradiction proves that our assumption
was wrong; qed.}.
\end{itemize}

Now, I claim that if $\sigma\in S_{n}$ satisfies $\sigma\notin\left\{
\operatorname*{id},z\right\}  $, then%
\begin{equation}
\left(
\begin{array}
[c]{c}%
\text{there exists an }i\in\left[  n\right]  \text{ satisfying}\\
i\neq\sigma\left(  i\right)  \text{ and }i\not \equiv \sigma\left(  i\right)
+1\operatorname{mod}n
\end{array}
\right)  . \label{sol.det.2diags.exi}%
\end{equation}

[\textit{Proof of (\ref{sol.det.2diags.exi}):} Let $\sigma\in S_{n}$ be such
that $\sigma\notin\left\{  \operatorname*{id},z\right\}  $. Thus, $\sigma
\neq\operatorname*{id}$ and $\sigma\neq z$.

The map $\sigma$ is an element of $S_{n}$, thus a permutation of the set
$\left\{  1,2,\ldots,n\right\}  $ (since $S_{n}$ is the set of all
permutations of the set $\left\{  1,2,\ldots,n\right\}  $). In other words,
$\sigma$ is a permutation of the set $\left[  n\right]  $ (since $\left\{
1,2,\ldots,n\right\}  =\left[  n\right]  $). Thus, $\sigma$ is a bijection
$\left[  n\right]  \rightarrow\left[  n\right]  $. Hence, the map $\sigma$ is
bijective, and therefore injective and surjective.

There exists a $J\in\left[  n\right]  $ such that $\sigma\left(  J\right)
\neq J$\ \ \ \ \footnote{\textit{Proof.} Assume the contrary. Thus, there
exists no $J\in\left[  n\right]  $ such that $\sigma\left(  J\right)  \neq J$.
In other words, every $J\in\left[  n\right]  $ satisfies $\sigma\left(
J\right)  =J$. Thus, every $J\in\left[  n\right]  $ satisfies $\sigma\left(
J\right)  =J=\operatorname*{id}\left(  J\right)  $. In other words,
$\sigma=\operatorname*{id}$. This contradicts $\sigma\neq\operatorname*{id}$.
This contradiction proves that our assumption was wrong. Qed.}. Let $j$ be the
smallest such $J$. Thus, $j$ is an element of $\left[  n\right]  $ satisfying
$\sigma\left(  j\right)  \neq j$\ \ \ \ \footnote{\textit{Proof.} We know that
$j$ is the smallest $J\in\left[  n\right]  $ such that $\sigma\left(
J\right)  \neq J$ (by the definition of $j$). Hence, $j$ is an element $J$ of
$\left[  n\right]  $ such that $\sigma\left(  J\right)  \neq J$. In other
words, $j$ is an element of $\left[  n\right]  $ satisfying $\sigma\left(
j\right)  \neq j$. Qed.}. Moreover,
\begin{equation}
\text{every }J\in\left[  n\right]  \text{ satisfying }\sigma\left(  J\right)
\neq J\text{ must satisfy }J\geq j \label{sol.det.2diags.exi.pf.min}%
\end{equation}
\ \ \ \ \footnote{\textit{Proof of (\ref{sol.det.2diags.exi.pf.min}):} We know
that $j$ is the \textbf{smallest} $J\in\left[  n\right]  $ such that
$\sigma\left(  J\right)  \neq J$ (by the definition of $j$). Hence, no
$J\in\left[  n\right]  $ such that $\sigma\left(  J\right)  \neq J$ can be
smaller than $j$. In other words, every $J\in\left[  n\right]  $ such that
$\sigma\left(  J\right)  \neq J$ must be $\geq j$. In other words, every
$J\in\left[  n\right]  $ satisfying $\sigma\left(  J\right)  \neq J$ must
satisfy $J\geq j$. This proves (\ref{sol.det.2diags.exi.pf.min}).}.

If $j\not \equiv \sigma\left(  j\right)  +1\operatorname{mod}n$, then
(\ref{sol.det.2diags.exi}) holds\footnote{\textit{Proof.} Assume the contrary.
Thus, $j\not \equiv \sigma\left(  j\right)  +1\operatorname{mod}n$. So we know
that $j\neq\sigma\left(  j\right)  $ (since $\sigma\left(  j\right)  \neq j$)
and $j\not \equiv \sigma\left(  j\right)  +1\operatorname{mod}n$.
Consequently, there exists an $i\in\left[  n\right]  $ satisfying $i\neq
\sigma\left(  i\right)  $ and $i\not \equiv \sigma\left(  i\right)
+1\operatorname{mod}n$ (namely, $i=j$). In other words,
(\ref{sol.det.2diags.exi}) holds. Qed.}. Hence, for the rest of this proof of
(\ref{sol.det.2diags.exi}), we can WLOG assume that we don't have
$j\not \equiv \sigma\left(  j\right)  +1\operatorname{mod}n$. Assume this.

We have $j\equiv\sigma\left(  j\right)  +1\operatorname{mod}n$ (since we don't
have $j\not \equiv \sigma\left(  j\right)  +1\operatorname{mod}n$). In other
words, $\sigma\left(  j\right)  \equiv j-1\operatorname{mod}n$.

We have $j=1$\ \ \ \ \footnote{\textit{Proof.} Assume the contrary. Thus,
$j\neq1$. Combining this with $j\in\left[  n\right]  =\left\{  1,2,\ldots
,n\right\}  $, we obtain $j\in\left\{  1,2,\ldots,n\right\}  \setminus\left\{
1\right\}  =\left\{  2,3,\ldots,n\right\}  $. Hence, $j-1\in\left\{
1,2,\ldots,n-1\right\}  \subseteq\left\{  1,2,\ldots,n\right\}  =\left[
n\right]  $. Also, $j-1<j$.
\par
If we had $\sigma\left(  j-1\right)  \neq j-1$, then we would obtain $j-1\geq
j$ (by (\ref{sol.det.2diags.exi.pf.min}), applied to $J=j-1$); but this would
contradict $j-1<j$. Hence, we cannot have $\sigma\left(  j-1\right)  \neq
j-1$. In other words, we must have $\sigma\left(  j-1\right)  =j-1$. Thus,
$\sigma\left(  j-1\right)  =j-1\equiv\sigma\left(  j\right)
\operatorname{mod}n$ (since $\sigma\left(  j\right)  \equiv
j-1\operatorname{mod}n$). Since both $\sigma\left(  j-1\right)  $ and
$\sigma\left(  j\right)  $ are elements of $\left[  n\right]  $, we can now
apply (\ref{sol.det.2diags.posrem.3}) to $\sigma\left(  j-1\right)  $ and
$\sigma\left(  j\right)  $ instead of $a$ and $b$. As a result, we obtain
$\sigma\left(  j-1\right)  =\sigma\left(  j\right)  $. Since the map $\sigma$
is injective, this shows that $j-1=j$. But this contradicts $j-1\neq j$. This
contradiction shows that our assumption was wrong, qed.}. Thus, $\sigma\left(
\underbrace{1}_{=j}\right)  =\sigma\left(  j\right)  \equiv\underbrace{j}%
_{=1}-1=0\equiv n\operatorname{mod}n$. Since both $\sigma\left(  1\right)  $
and $n$ belong to $\left\{  1,2,\ldots,n\right\}  =\left[  n\right]  $, we can
now apply (\ref{sol.det.2diags.posrem.3}) to $\sigma\left(  1\right)  $ and
$n$ instead of $a$ and $b$. As a result, we obtain $\sigma\left(  1\right)
=n$.

There exists a $p\in\left[  n\right]  $ such that $\sigma\left(  p\right)
\neq z\left(  p\right)  $\ \ \ \ \footnote{\textit{Proof.} Assume the
contrary. Thus, there exists no $p\in\left[  n\right]  $ such that
$\sigma\left(  p\right)  \neq z\left(  p\right)  $. In other words, every
$p\in\left[  n\right]  $ satisfies $\sigma\left(  p\right)  =z\left(
p\right)  $. In other words, $\sigma=z$. This contradicts $\sigma\neq z$. This
contradiction shows that our assumption was wrong. Qed.}. Consider this $p$.

We have $p\not \equiv \sigma\left(  p\right)  +1\operatorname{mod}%
n$\ \ \ \ \footnote{\textit{Proof.} Assume the contrary. Thus, $p\equiv
\sigma\left(  p\right)  +1\operatorname{mod}n$. Hence, $\sigma\left(
p\right)  \equiv p-1\operatorname{mod}n$. But (\ref{sol.det.2diags.zmod})
(applied to $i=p$) yields $z\left(  p\right)  \equiv p-1\operatorname{mod}n$.
Hence, $\sigma\left(  p\right)  \equiv p-1\equiv z\left(  p\right)
\operatorname{mod}n$. Since both $\sigma$ and $z$ are elements of $S_{n}$, we
have $\sigma\left(  p\right)  \in\left[  n\right]  $ and $z\left(  p\right)
\in\left[  n\right]  $. Thus, (\ref{sol.det.2diags.posrem.3}) (applied to
$\sigma\left(  p\right)  $ and $z\left(  p\right)  $ instead of $a$ and $b$)
shows that $\sigma\left(  p\right)  =z\left(  p\right)  $ (since
$\sigma\left(  p\right)  \equiv z\left(  p\right)  \operatorname{mod}n$). This
contradicts $\sigma\left(  p\right)  \neq z\left(  p\right)  $. This
contradiction proves that our assumption was wrong, qed.}. If we have
$p\neq\sigma\left(  p\right)  $, then (\ref{sol.det.2diags.exi})
holds\footnote{\textit{Proof.} Assume that $p\neq\sigma\left(  p\right)  $. So
we know that $p\neq\sigma\left(  p\right)  $ and $p\not \equiv \sigma\left(
p\right)  +1\operatorname{mod}n$. Consequently, there exists an $i\in\left[
n\right]  $ satisfying $i\neq\sigma\left(  i\right)  $ and $i\not \equiv
\sigma\left(  i\right)  +1\operatorname{mod}n$ (namely, $i=p$). In other
words, (\ref{sol.det.2diags.exi}) holds. Qed.}. Hence, for the rest of this
proof of (\ref{sol.det.2diags.exi}), we can WLOG assume that we don't have
$p\neq\sigma\left(  p\right)  $. Assume this.

We have $p=\sigma\left(  p\right)  $ (since we don't have $p\neq\sigma\left(
p\right)  $). In other words, $\sigma\left(  p\right)  =p$. Hence, there
exists a $K\in\left[  n\right]  $ such that $\sigma\left(  K\right)  =K$
(namely, $K=p$). Let $k$ be the largest such $K$. Thus, $k$ is an element of
$\left[  n\right]  $ satisfying $\sigma\left(  k\right)  =k$%
\ \ \ \ \footnote{\textit{Proof.} We know that $k$ is the largest $K\in\left[
n\right]  $ such that $\sigma\left(  K\right)  =K$ (by the definition of $k$).
Hence, $k$ is an element $K$ of $\left[  n\right]  $ such that $\sigma\left(
K\right)  =K$. In other words, $k$ is an element of $\left[  n\right]  $
satisfying $\sigma\left(  k\right)  =k$. Qed.}. Moreover,
\begin{equation}
\text{every }K\in\left[  n\right]  \text{ satisfying }\sigma\left(  K\right)
=K\text{ must satisfy }K\leq k \label{sol.det.2diags.exi.pf.max}%
\end{equation}
\ \ \ \ \footnote{\textit{Proof of (\ref{sol.det.2diags.exi.pf.max}):} We know
that $k$ is the \textbf{largest} $K\in\left[  n\right]  $ such that
$\sigma\left(  K\right)  =K$ (by the definition of $k$). Hence, no
$K\in\left[  n\right]  $ such that $\sigma\left(  K\right)  =K$ can be larger
than $k$. In other words, every $K\in\left[  n\right]  $ such that
$\sigma\left(  K\right)  =K$ must be $\leq k$. In other words, every
$K\in\left[  n\right]  $ satisfying $\sigma\left(  K\right)  =K$ must satisfy
$K\leq k$. This proves (\ref{sol.det.2diags.exi.pf.max}).}.

We have $n\neq1$ and thus $\sigma\left(  n\right)  \neq\sigma\left(  1\right)
$ (since the map $\sigma$ is injective). Thus, $\sigma\left(  n\right)
\neq\sigma\left(  1\right)  =n$.

If we had $k=n$, then we would have $\sigma\left(  \underbrace{n}_{=k}\right)
=\sigma\left(  k\right)  =k=n$, which would contradict $\sigma\left(
n\right)  \neq n$. Thus, we cannot have $k=n$. Hence, we have $k\neq n$. Thus,
$k+1\in\left[  n\right]  $\ \ \ \ \footnote{\textit{Proof.} Combining
$k\in\left[  n\right]  =\left\{  1,2,\ldots,n\right\}  $ with $k\neq n$, we
obtain $k\in\left\{  1,2,\ldots,n\right\}  \setminus\left\{  n\right\}
=\left\{  1,2,\ldots,n-1\right\}  $, so that $k+1\in\left\{  2,3,\ldots
,n\right\}  \subseteq\left\{  1,2,\ldots,n\right\}  =\left[  n\right]  $,
qed.}. Hence, $\sigma\left(  k+1\right)  $ is well-defined. If we had
$\sigma\left(  k+1\right)  =k+1$, then we would obtain $k+1\leq k$ (by
(\ref{sol.det.2diags.exi.pf.max}), applied to $K=k+1$), which would contradict
$k+1>k$. Hence, we cannot have $\sigma\left(  k+1\right)  =k+1$. We therefore
must have $\sigma\left(  k+1\right)  \neq k+1$. In other words, $k+1\neq
\sigma\left(  k+1\right)  $.

But we also have $k+1\not \equiv \sigma\left(  k+1\right)
+1\operatorname{mod}n$\ \ \ \ \footnote{\textit{Proof.} Assume the contrary.
Thus, $k+1\equiv\sigma\left(  k+1\right)  +1\operatorname{mod}n$. Subtracting
$1$ from both sides of this congruence, we obtain $k\equiv\sigma\left(
k+1\right)  \operatorname{mod}n$. Since both $k$ and $\sigma\left(
k+1\right)  $ belong to $\left[  n\right]  $, we can therefore apply
(\ref{sol.det.2diags.posrem.3}) to $k$ and $\sigma\left(  k+1\right)  $
instead of $a$ and $b$. As a result, we obtain $k=\sigma\left(  k+1\right)  $.
Hence, $\sigma\left(  k+1\right)  =k=\sigma\left(  k\right)  $ (since
$\sigma\left(  k\right)  =k$). Since the map $\sigma$ is injective, this
yields that $k+1=k$. But this contradicts $k+1>k$. This contradiction proves
that our assumption was false, qed.}. Hence, there exists an $i\in\left[
n\right]  $ satisfying $i\neq\sigma\left(  i\right)  $ and $i\not \equiv
\sigma\left(  i\right)  +1\operatorname{mod}n$ (namely, $i=k+1$). In other
words, (\ref{sol.det.2diags.exi}) holds. This completes the proof of
(\ref{sol.det.2diags.exi}).]

Let us now write our matrix $A$ in the form $A=\left(  a_{i,j}\right)  _{1\leq
i\leq n,\ 1\leq j\leq n}$. Then,%
\[
\left(  a_{i,j}\right)  _{1\leq i\leq n,\ 1\leq j\leq n}=A=\left(
\begin{cases}
a_{j}, & \text{if }i=j;\\
b_{j}, & \text{if }i\equiv j+1\operatorname{mod}n;\\
0, & \text{otherwise}%
\end{cases}
\right)  _{1\leq i\leq n,\ 1\leq j\leq n}.
\]
In other words, we have%
\begin{equation}
a_{i,j}=%
\begin{cases}
a_{j}, & \text{if }i=j;\\
b_{j}, & \text{if }i\equiv j+1\operatorname{mod}n;\\
0, & \text{otherwise}%
\end{cases}
\label{sol.det.2diags.aij}%
\end{equation}
for every $\left(  i,j\right)  \in\left\{  1,2,\ldots,n\right\}  ^{2}$.

It is now easy to see that if $\sigma\in S_{n}$ satisfies $\sigma
\notin\left\{  \operatorname*{id},z\right\}  $, then
\begin{equation}
\prod_{i=1}^{n}a_{i,\sigma\left(  i\right)  }=0 \label{sol.det.2diags.0}%
\end{equation}
\footnote{\textit{Proof of (\ref{sol.det.2diags.0}):} Let $\sigma\in S_{n}$ be
such that $\sigma\notin\left\{  \operatorname*{id},z\right\}  $. Thus, there
exists an $i\in\left[  n\right]  $ satisfying $i\neq\sigma\left(  i\right)  $
and $i\not \equiv \sigma\left(  i\right)  +1\operatorname{mod}n$ (because of
(\ref{sol.det.2diags.exi})). Consider this $i$. We have neither $i=\sigma
\left(  i\right)  $ nor $i\equiv\sigma\left(  i\right)  +1\operatorname{mod}n$
(since we have $i\neq\sigma\left(  i\right)  $ and $i\not \equiv \sigma\left(
i\right)  +1\operatorname{mod}n$).
\par
From (\ref{sol.det.2diags.aij}) (applied to $\left(  i,\sigma\left(  i\right)
\right)  $ instead of $\left(  i,j\right)  $), we obtain%
\[
a_{i,\sigma\left(  i\right)  }=%
\begin{cases}
a_{\sigma\left(  i\right)  }, & \text{if }i=\sigma\left(  i\right)  ;\\
b_{\sigma\left(  i\right)  }, & \text{if }i\equiv\sigma\left(  i\right)
+1\operatorname{mod}n;\\
0, & \text{otherwise}%
\end{cases}
=0
\]
(since we have neither $i=\sigma\left(  i\right)  $ nor $i\equiv\sigma\left(
i\right)  +1\operatorname{mod}n$).
\par
Now, let us forget that we fixed $i$. We thus have found an $i\in\left[
n\right]  $ satisfying $a_{i,\sigma\left(  i\right)  }=0$. Thus, there exists
an $i\in\left[  n\right]  $ satisfying $a_{i,\sigma\left(  i\right)  }=0$. In
other words, at least one factor of the product $\prod_{i\in\left[  n\right]
}a_{i,\sigma\left(  i\right)  }$ equals $0$. Therefore, the whole product
$\prod_{i\in\left[  n\right]  }a_{i,\sigma\left(  i\right)  }$ equals $0$
(because if at least one factor of a product equals $0$, then the whole
product equals $0$). In other words, $\prod_{i\in\left[  n\right]
}a_{i,\sigma\left(  i\right)  }=0$.
\par
Now,
\[
\underbrace{\prod_{i=1}^{n}}_{\substack{=\prod_{i\in\left\{  1,2,\ldots
,n\right\}  }=\prod_{i\in\left[  n\right]  }\\\text{(since }\left\{
1,2,\ldots,n\right\}  =\left[  n\right]  \text{)}}}a_{i,\sigma\left(
i\right)  }=\prod_{i\in\left[  n\right]  }a_{i,\sigma\left(  i\right)  }=0.
\]
This proves (\ref{sol.det.2diags.0}).}.

Let us also notice that%
\begin{equation}
\prod_{i=1}^{n}a_{i,z\left(  i\right)  }=\prod_{i=1}^{n}b_{i}
\label{sol.det.2diags.1}%
\end{equation}
\footnote{\textit{Proof of (\ref{sol.det.2diags.1}):} The map $z$ is a
permutation of $\left[  n\right]  $. In other words, the map $z$ is a
permutation of $\left\{  1,2,\ldots,n\right\}  $ (since $\left[  n\right]
=\left\{  1,2,\ldots,n\right\}  $), thus a bijection $\left\{  1,2,\ldots
,n\right\}  \rightarrow\left\{  1,2,\ldots,n\right\}  $. Hence, we can
substitute $i$ for $z\left(  i\right)  $ in the product $\prod_{i\in\left\{
1,2,\ldots,n\right\}  }b_{z\left(  i\right)  }$. We thus obtain $\prod
_{i\in\left\{  1,2,\ldots,n\right\}  }b_{z\left(  i\right)  }=\prod
_{i\in\left\{  1,2,\ldots,n\right\}  }b_{i}$.
\par
But for every $i\in\left\{  1,2,\ldots,n\right\}  $, we have%
\begin{align*}
a_{i,z\left(  i\right)  }  &  =%
\begin{cases}
a_{z\left(  i\right)  }, & \text{if }i=z\left(  i\right)  ;\\
b_{z\left(  i\right)  }, & \text{if }i\equiv z\left(  i\right)
+1\operatorname{mod}n;\\
0, & \text{otherwise}%
\end{cases}
\\
&  \ \ \ \ \ \ \ \ \ \ \left(  \text{by (\ref{sol.det.2diags.aij}), applied to
}\left(  i,z\left(  i\right)  \right)  \text{ instead of }\left(  i,j\right)
\right) \\
&  =b_{z\left(  i\right)  }\ \ \ \ \ \ \ \ \ \ \left(
\begin{array}
[c]{c}%
\text{since we don't have }i=z\left(  i\right)  \text{ (because of
(\ref{sol.det.2diags.zmod0})),}\\
\text{but we do have }i\equiv z\left(  i\right)  +1\operatorname{mod}n\text{
(by (\ref{sol.det.2diags.zmod}))}%
\end{array}
\right)  .
\end{align*}
Thus,%
\[
\underbrace{\prod_{i=1}^{n}}_{=\prod_{i\in\left\{  1,2,\ldots,n\right\}  }%
}\underbrace{a_{i,z\left(  i\right)  }}_{=b_{z\left(  i\right)  }}=\prod
_{i\in\left\{  1,2,\ldots,n\right\}  }b_{z\left(  i\right)  }%
=\underbrace{\prod_{i\in\left\{  1,2,\ldots,n\right\}  }}_{=\prod_{i=1}^{n}%
}b_{i}=\prod_{i=1}^{n}b_{i}.
\]
This proves (\ref{sol.det.2diags.1}).}.

Now, (\ref{eq.det.eq.2}) becomes%
\begin{align*}
\det A  &  =\sum_{\sigma\in S_{n}}\left(  -1\right)  ^{\sigma}\prod_{i=1}%
^{n}a_{i,\sigma\left(  i\right)  }\\
&  =\underbrace{\left(  -1\right)  ^{\operatorname*{id}}}_{=1}\prod_{i=1}%
^{n}\underbrace{a_{i,\operatorname*{id}\left(  i\right)  }}%
_{\substack{=a_{i,i}\\\text{(since }\operatorname*{id}\left(  i\right)
=i\text{)}}}+\underbrace{\left(  -1\right)  ^{z}}_{\substack{=\left(
-1\right)  ^{n-1}\\\text{(by (\ref{sol.det.2diags.signz}))}}}\prod_{i=1}%
^{n}a_{i,z\left(  i\right)  }+\sum_{\substack{\sigma\in S_{n};\\\sigma
\notin\left\{  \operatorname*{id},z\right\}  }}\left(  -1\right)  ^{\sigma
}\underbrace{\prod_{i=1}^{n}a_{i,\sigma\left(  i\right)  }}%
_{\substack{=0\\\text{(by (\ref{sol.det.2diags.0}))}}}\\
&  \ \ \ \ \ \ \ \ \ \ \left(
\begin{array}
[c]{c}%
\text{here, we have split off the addends for }\sigma=\operatorname*{id}\text{
and}\\
\text{for }\sigma=z\text{ from the sum (since }z\neq\operatorname*{id}\text{)}%
\end{array}
\right) \\
&  =\prod_{i=1}^{n}a_{i,i}+\left(  -1\right)  ^{n-1}\prod_{i=1}^{n}%
a_{i,z\left(  i\right)  }+\underbrace{\sum_{\substack{\sigma\in S_{n}%
;\\\sigma\notin\left\{  \operatorname*{id},z\right\}  }}\left(  -1\right)
^{\sigma}0}_{=0}\\
&  =\prod_{i=1}^{n}\underbrace{a_{i,i}}_{\substack{=%
\begin{cases}
a_{i}, & \text{if }i=i;\\
b_{i}, & \text{if }i\equiv i+1\operatorname{mod}n;\\
0, & \text{otherwise}%
\end{cases}
\\\text{(by (\ref{sol.det.2diags.aij}), applied to }\left(  i,i\right)  \text{
instead of }\left(  i,j\right)  \text{)}}}+\left(  -1\right)  ^{n-1}%
\underbrace{\prod_{i=1}^{n}a_{i,z\left(  i\right)  }}_{\substack{=\prod
_{i=1}^{n}b_{i}\\\text{(by (\ref{sol.det.2diags.1}))}}}\\
&  =\prod_{i=1}^{n}\underbrace{%
\begin{cases}
a_{i}, & \text{if }i=i;\\
b_{i}, & \text{if }i\equiv i+1\operatorname{mod}n;\\
0, & \text{otherwise}%
\end{cases}
}_{\substack{=a_{i}\\\text{(since }i=i\text{)}}}+\left(  -1\right)
^{n-1}\prod_{i=1}^{n}b_{i}\\
&  =\underbrace{\prod_{i=1}^{n}a_{i}}_{=a_{1}a_{2}\cdots a_{n}}+\left(
-1\right)  ^{n-1}\underbrace{\prod_{i=1}^{n}b_{i}}_{=b_{1}b_{2}\cdots b_{n}}\\
&  =a_{1}a_{2}\cdots a_{n}+\left(  -1\right)  ^{n-1}b_{1}b_{2}\cdots b_{n}.
\end{align*}
This solves Exercise \ref{exe.det.2diags}.
\end{proof}
\end{verlong}

\subsection{Solution to Exercise \ref{exe.altern.STAS}}

Before we come to the solution of Exercise \ref{exe.altern.STAS}, let us first
state some simple lemmas:

\begin{lemma}
\label{lem.sol.altern.STAS.1}Let $n\in\mathbb{N}$, $m\in\mathbb{N}$,
$p\in\mathbb{N}$ and $q\in\mathbb{N}$. Let $A=\left(  a_{i,j}\right)  _{1\leq
i\leq n,\ 1\leq j\leq m}\in\mathbb{K}^{n\times m}$, $B=\left(  b_{i,j}\right)
_{1\leq i\leq m,\ 1\leq j\leq p}\in\mathbb{K}^{m\times p}$ and $C=\left(
c_{i,j}\right)  _{1\leq i\leq p,\ 1\leq j\leq q}\in\mathbb{K}^{p\times q}$.
Then,%
\[
ABC=\left(  \sum_{k=1}^{m}\sum_{\ell=1}^{p}a_{i,k}b_{k,\ell}c_{\ell,j}\right)
_{1\leq i\leq n,\ 1\leq j\leq q}.
\]

\end{lemma}

\begin{proof}
[Proof of Lemma \ref{lem.sol.altern.STAS.1}.]We have $B=\left(  b_{i,j}%
\right)  _{1\leq i\leq m,\ 1\leq j\leq p}$ and $C=\left(  c_{i,j}\right)
_{1\leq i\leq p,\ 1\leq j\leq q}$. Thus, the definition of the product $BC$
yields%
\[
BC=\left(  \underbrace{\sum_{k=1}^{p}b_{i,k}c_{k,j}}_{\substack{=\sum_{\ell
=1}^{p}b_{i,\ell}c_{\ell,j}\\\text{(here, we renamed the}\\\text{summation
index }k\text{ as }\ell\text{)}}}\right)  _{1\leq i\leq m,\ 1\leq j\leq
q}=\left(  \sum_{\ell=1}^{p}b_{i,\ell}c_{\ell,j}\right)  _{1\leq i\leq
m,\ 1\leq j\leq q}.
\]

Now, we have $A=\left(  a_{i,j}\right)  _{1\leq i\leq n,\ 1\leq j\leq m}$ and
$BC=\left(  \sum_{\ell=1}^{p}b_{i,\ell}c_{\ell,j}\right)  _{1\leq i\leq
m,\ 1\leq j\leq q}$. Hence, the definition of the product $A\left(  BC\right)
$ yields%
\[
A\left(  BC\right)  =\left(  \sum_{k=1}^{m}\underbrace{a_{i,k}\left(
\sum_{\ell=1}^{p}b_{k,\ell}c_{\ell,j}\right)  }_{=\sum_{\ell=1}^{p}%
a_{i,k}b_{k,\ell}c_{\ell,j}}\right)  _{1\leq i\leq n,\ 1\leq j\leq q}=\left(
\sum_{k=1}^{m}\sum_{\ell=1}^{p}a_{i,k}b_{k,\ell}c_{\ell,j}\right)  _{1\leq
i\leq n,\ 1\leq j\leq q}.
\]
Thus,%
\[
ABC=A\left(  BC\right)  =\left(  \sum_{k=1}^{m}\sum_{\ell=1}^{p}%
a_{i,k}b_{k,\ell}c_{\ell,j}\right)  _{1\leq i\leq n,\ 1\leq j\leq q}.
\]
This proves Lemma \ref{lem.sol.altern.STAS.1}.
\end{proof}

\begin{lemma}
\label{lem.sol.altern.STAS.2}Let $n\in\mathbb{N}$ and $m\in\mathbb{N}$. Let
$A=\left(  a_{i,j}\right)  _{1\leq i\leq n,\ 1\leq j\leq n}$ be an $n\times
n$-matrix. Let $S=\left(  s_{i,j}\right)  _{1\leq i\leq n,\ 1\leq j\leq m}$ be
an $n\times m$-matrix. Then,%
\[
S^{T}AS=\left(  \sum_{\left(  k,\ell\right)  \in\left\{  1,2,\ldots,n\right\}
^{2}}s_{k,i}s_{\ell,j}a_{k,\ell}\right)  _{1\leq i\leq m,\ 1\leq j\leq m}.
\]

\end{lemma}

\begin{proof}
[Proof of Lemma \ref{lem.sol.altern.STAS.2}.]We have $S=\left(  s_{i,j}%
\right)  _{1\leq i\leq n,\ 1\leq j\leq m}$. Thus, $S^{T}=\left(
s_{j,i}\right)  _{1\leq i\leq m,\ 1\leq j\leq n}$ (by the definition of
$S^{T}$). Recall also that $A=\left(  a_{i,j}\right)  _{1\leq i\leq n,\ 1\leq
j\leq n}$ and $S=\left(  s_{i,j}\right)  _{1\leq i\leq n,\ 1\leq j\leq m}$.
Hence, Lemma \ref{lem.sol.altern.STAS.1} (applied to $m$, $n$, $n$, $m$,
$S^{T}$, $A$, $S$, $s_{j,i}$, $a_{i,j}$ and $s_{i,j}$ instead of $n$, $m$,
$p$, $q$, $A$, $B$, $C$, $a_{i,j}$, $b_{i,j}$ and $c_{i,j}$) yields%
\begin{align*}
S^{T}AS  &  =\left(  \underbrace{\sum_{k=1}^{n}}_{=\sum_{k\in\left\{
1,2,\ldots,n\right\}  }}\ \ \underbrace{\sum_{\ell=1}^{n}}_{=\sum_{\ell
\in\left\{  1,2,\ldots,n\right\}  }}s_{k,i}\underbrace{a_{k,\ell}s_{\ell,j}%
}_{=s_{\ell,j}a_{k,\ell}}\right)  _{1\leq i\leq m,\ 1\leq j\leq m}\\
&  =\left(  \underbrace{\sum_{k\in\left\{  1,2,\ldots,n\right\}  }%
\ \ \sum_{\ell\in\left\{  1,2,\ldots,n\right\}  }}_{=\sum_{\left(
k,\ell\right)  \in\left\{  1,2,\ldots,n\right\}  ^{2}}}s_{k,i}s_{\ell
,j}a_{k,\ell}\right)  _{1\leq i\leq m,\ 1\leq j\leq m}\\
&  =\left(  \sum_{\left(  k,\ell\right)  \in\left\{  1,2,\ldots,n\right\}
^{2}}s_{k,i}s_{\ell,j}a_{k,\ell}\right)  _{1\leq i\leq m,\ 1\leq j\leq m}.
\end{align*}
This proves Lemma \ref{lem.sol.altern.STAS.2}.
\end{proof}

\begin{lemma}
\label{lem.sol.altern.STAS.3}Let $n\in\mathbb{N}$. Let $A=\left(
a_{i,j}\right)  _{1\leq i\leq n,\ 1\leq j\leq n}$ be an alternating $n\times
n$-matrix.

\textbf{(a)} Every $i\in\left\{  1,2,\ldots,n\right\}  $ satisfies $a_{i,i}=0$.

\textbf{(b)} Every $\left(  i,j\right)  \in\left\{  1,2,\ldots,n\right\}
^{2}$ satisfy $a_{i,j}=-a_{j,i}$.
\end{lemma}

\begin{proof}
[Proof of Lemma \ref{lem.sol.altern.STAS.3}.]Recall that the matrix $A$ is
alternating if and only if it satisfies $A^{T}=-A$ and $\left(  a_{i,i}%
=0\text{ for all }i\in\left\{  1,2,\ldots,n\right\}  \right)  $ (by the
definition of \textquotedblleft alternating\textquotedblright). Hence, the
matrix $A$ satisfies $A^{T}=-A$ and $\left(  a_{i,i}=0\text{ for all }%
i\in\left\{  1,2,\ldots,n\right\}  \right)  $ (since $A$ is alternating).

We have $\left(  a_{i,i}=0\text{ for all }i\in\left\{  1,2,\ldots,n\right\}
\right)  $. In other words, every $i\in\left\{  1,2,\ldots,n\right\}  $
satisfies $a_{i,i}=0$. This proves Lemma \ref{lem.sol.altern.STAS.3}
\textbf{(a)}.

\textbf{(b)} We have $A=\left(  a_{i,j}\right)  _{1\leq i\leq n,\ 1\leq j\leq
n}$, and thus $A^{T}=\left(  a_{j,i}\right)  _{1\leq i\leq n,\ 1\leq j\leq n}$
(by the definition of $A^{T}$). Hence,%
\[
\left(  a_{j,i}\right)  _{1\leq i\leq n,\ 1\leq j\leq n}=A^{T}=-\underbrace{A}%
_{=\left(  a_{i,j}\right)  _{1\leq i\leq n,\ 1\leq j\leq n}}=-\left(
a_{i,j}\right)  _{1\leq i\leq n,\ 1\leq j\leq n}=\left(  -a_{i,j}\right)
_{1\leq i\leq n,\ 1\leq j\leq n}.
\]
In other words,%
\begin{equation}
a_{j,i}=-a_{i,j}\ \ \ \ \ \ \ \ \ \ \text{for every }\left(  i,j\right)
\in\left\{  1,2,\ldots,n\right\}  ^{2}. \label{pf.lem.sol.altern.STAS.3.1}%
\end{equation}

Now, let $\left(  i,j\right)  \in\left\{  1,2,\ldots,n\right\}  ^{2}$. Hence,
(\ref{pf.lem.sol.altern.STAS.3.1}) yields $a_{j,i}=-a_{i,j}$. In other words,
$a_{i,j}=-a_{j,i}$. This proves Lemma \ref{lem.sol.altern.STAS.3} \textbf{(b)}.
\end{proof}

\begin{proof}
[Solution to Exercise \ref{exe.altern.STAS}.]Let $B$ be the $m\times m$-matrix
$S^{T}AS$. Thus, $B=S^{T}AS$. Write the $m\times m$-matrix $B$ in the form
$B=\left(  b_{i,j}\right)  _{1\leq i\leq m,\ 1\leq j\leq m}$.

Write the $n\times n$-matrix $A$ in the form $A=\left(  a_{i,j}\right)
_{1\leq i\leq n,\ 1\leq j\leq n}$. Write the $n\times m$-matrix $S$ in the
form $S=\left(  s_{i,j}\right)  _{1\leq i\leq n,\ 1\leq j\leq m}$. Then,%
\begin{equation}
B=S^{T}AS=\left(  \sum_{\left(  k,\ell\right)  \in\left\{  1,2,\ldots
,n\right\}  ^{2}}s_{k,i}s_{\ell,j}a_{k,\ell}\right)  _{1\leq i\leq m,\ 1\leq
j\leq m} \label{sol.altern.STAS.1}%
\end{equation}
(by Lemma \ref{lem.sol.altern.STAS.2}).

But $B=\left(  b_{i,j}\right)  _{1\leq i\leq m,\ 1\leq j\leq m}$, so that%
\[
\left(  b_{i,j}\right)  _{1\leq i\leq m,\ 1\leq j\leq m}=B=\left(
\sum_{\left(  k,\ell\right)  \in\left\{  1,2,\ldots,n\right\}  ^{2}}%
s_{k,i}s_{\ell,j}a_{k,\ell}\right)  _{1\leq i\leq m,\ 1\leq j\leq m}%
\]
(by (\ref{sol.altern.STAS.1})). In other words,%
\begin{equation}
b_{i,j}=\sum_{\left(  k,\ell\right)  \in\left\{  1,2,\ldots,n\right\}  ^{2}%
}s_{k,i}s_{\ell,j}a_{k,\ell}\ \ \ \ \ \ \ \ \ \ \text{for every }\left(
i,j\right)  \in\left\{  1,2,\ldots,m\right\}  ^{2}. \label{sol.altern.STAS.4}%
\end{equation}

\begin{vershort}
But
\[
\left(  b_{j,i}\right)  _{1\leq i\leq m,\ 1\leq j\leq m}=\left(
-b_{i,j}\right)  _{1\leq i\leq m,\ 1\leq j\leq m}%
\]
\footnote{\textit{Proof:} Every $\left(  i,j\right)  \in\left\{
1,2,\ldots,m\right\}  ^{2}$ satisfies
\begin{align*}
b_{j,i}  &  =\sum_{\left(  k,\ell\right)  \in\left\{  1,2,\ldots,n\right\}
^{2}}s_{k,j}s_{\ell,i}a_{k,\ell}\ \ \ \ \ \ \ \ \ \ \left(  \text{by
(\ref{sol.altern.STAS.4}) (applied to }\left(  j,i\right)  \text{ instead of
}\left(  i,j\right)  \text{)}\right) \\
&  =\underbrace{\sum_{\left(  \ell,k\right)  \in\left\{  1,2,\ldots,n\right\}
^{2}}}_{\substack{=\sum_{\left(  k,\ell\right)  \in\left\{  1,2,\ldots
,n\right\}  ^{2}}}}\underbrace{s_{\ell,j}s_{k,i}}_{=s_{k,i}s_{\ell,j}%
}\underbrace{a_{\ell,k}}_{\substack{=-a_{k,\ell}\\\text{(by Lemma
\ref{lem.sol.altern.STAS.3} \textbf{(b)}}\\\text{(applied to }\left(
\ell,k\right)  \\\text{instead of }\left(  i,j\right)  \text{))}%
}}\ \ \ \ \ \ \ \ \ \ \left(
\begin{array}
[c]{c}%
\text{here, we have renamed the}\\
\text{summation index }\left(  k,\ell\right)  \text{ as }\left(
\ell,k\right)
\end{array}
\right) \\
&  =\sum_{\left(  k,\ell\right)  \in\left\{  1,2,\ldots,n\right\}  ^{2}%
}s_{k,i}s_{\ell,j}\left(  -a_{k,\ell}\right)  =-\underbrace{\sum_{\left(
k,\ell\right)  \in\left\{  1,2,\ldots,n\right\}  ^{2}}s_{k,i}s_{\ell
,j}a_{k,\ell}}_{\substack{=b_{i,j}\\\text{(by (\ref{sol.altern.STAS.4}))}%
}}=-b_{i,j}.
\end{align*}
In other words, $\left(  b_{j,i}\right)  _{1\leq i\leq m,\ 1\leq j\leq
m}=\left(  -b_{i,j}\right)  _{1\leq i\leq m,\ 1\leq j\leq m}$.}. But
$B=\left(  b_{i,j}\right)  _{1\leq i\leq m,\ 1\leq j\leq m}$, and thus
$B^{T}=\left(  b_{j,i}\right)  _{1\leq i\leq m,\ 1\leq j\leq m}$ (by the
definition of $B^{T}$). Hence,%
\[
B^{T}=\left(  b_{j,i}\right)  _{1\leq i\leq m,\ 1\leq j\leq m}=\left(
-b_{i,j}\right)  _{1\leq i\leq m,\ 1\leq j\leq m}.
\]

\end{vershort}

\begin{verlong}
But every $\left(  i,j\right)  \in\left\{  1,2,\ldots,m\right\}  ^{2}$
satisfies%
\begin{equation}
b_{j,i}=-b_{i,j} \label{sol.altern.STAS.5}%
\end{equation}
\footnote{\textit{Proof of (\ref{sol.altern.STAS.5}):} Let $\left(
i,j\right)  \in\left\{  1,2,\ldots,m\right\}  ^{2}$. Thus, $i\in\left\{
1,2,\ldots,m\right\}  $ and $j\in\left\{  1,2,\ldots,m\right\}  $. Hence,
$j\in\left\{  1,2,\ldots,m\right\}  $ and $i\in\left\{  1,2,\ldots,m\right\}
$. In other words, $\left(  j,i\right)  \in\left\{  1,2,\ldots,m\right\}
^{2}$. Now, (\ref{sol.altern.STAS.4}) (applied to $\left(  j,i\right)  $
instead of $\left(  i,j\right)  $) yields%
\begin{align*}
b_{j,i}  &  =\sum_{\left(  k,\ell\right)  \in\left\{  1,2,\ldots,n\right\}
^{2}}s_{k,j}s_{\ell,i}a_{k,\ell}\\
&  =\underbrace{\sum_{\left(  \ell,k\right)  \in\left\{  1,2,\ldots,n\right\}
^{2}}}_{\substack{=\sum_{\ell\in\left\{  1,2,\ldots,n\right\}  }\sum
_{k\in\left\{  1,2,\ldots,n\right\}  }\\=\sum_{k\in\left\{  1,2,\ldots
,n\right\}  }\sum_{\ell\in\left\{  1,2,\ldots,n\right\}  }\\=\sum_{\left(
k,\ell\right)  \in\left\{  1,2,\ldots,n\right\}  ^{2}}}}\underbrace{s_{\ell
,j}s_{k,i}}_{=s_{k,i}s_{\ell,j}}\underbrace{a_{\ell,k}}_{\substack{=-a_{k,\ell
}\\\text{(by Lemma \ref{lem.sol.altern.STAS.3} \textbf{(b)}}\\\text{(applied
to }\left(  \ell,k\right)  \\\text{instead of }\left(  i,j\right)  \text{))}%
}}\ \ \ \ \ \ \ \ \ \ \left(
\begin{array}
[c]{c}%
\text{here, we have renamed the}\\
\text{summation index }\left(  k,\ell\right)  \text{ as }\left(
\ell,k\right)
\end{array}
\right) \\
&  =\sum_{\left(  k,\ell\right)  \in\left\{  1,2,\ldots,n\right\}  ^{2}%
}s_{k,i}s_{\ell,j}\left(  -a_{k,\ell}\right)  =-\underbrace{\sum_{\left(
k,\ell\right)  \in\left\{  1,2,\ldots,n\right\}  ^{2}}s_{k,i}s_{\ell
,j}a_{k,\ell}}_{\substack{=b_{i,j}\\\text{(by (\ref{sol.altern.STAS.4}))}%
}}=-b_{i,j}.
\end{align*}
This proves (\ref{sol.altern.STAS.5}).}. But $B=\left(  b_{i,j}\right)
_{1\leq i\leq m,\ 1\leq j\leq m}$, and thus $B^{T}=\left(  b_{j,i}\right)
_{1\leq i\leq m,\ 1\leq j\leq m}$ (by the definition of $B^{T}$). Hence,%
\[
B^{T}=\left(  \underbrace{b_{j,i}}_{=-b_{i,j}}\right)  _{1\leq i\leq m,\ 1\leq
j\leq m}=\left(  -b_{i,j}\right)  _{1\leq i\leq m,\ 1\leq j\leq m}.
\]

\end{verlong}

Comparing this with%
\[
-\underbrace{B}_{=\left(  b_{i,j}\right)  _{1\leq i\leq m,\ 1\leq j\leq m}%
}=-\left(  b_{i,j}\right)  _{1\leq i\leq m,\ 1\leq j\leq m}=\left(
-b_{i,j}\right)  _{1\leq i\leq m,\ 1\leq j\leq m},
\]
we obtain $B^{T}=-B$.

\begin{vershort}
Let $i\in\left\{  1,2,\ldots,m\right\}  $. Then, (\ref{sol.altern.STAS.4})
(applied to $\left(  i,i\right)  $ instead of $\left(  i,j\right)  $) yields%
\begin{align*}
b_{i,i}  &  =\sum_{\left(  k,\ell\right)  \in\left\{  1,2,\ldots,n\right\}
^{2}}s_{k,i}s_{\ell,i}a_{k,\ell}\\
&  =\sum_{\substack{\left(  k,\ell\right)  \in\left\{  1,2,\ldots,n\right\}
^{2};\\k<\ell}}s_{k,i}s_{\ell,i}a_{k,\ell}+\sum_{\substack{\left(
k,\ell\right)  \in\left\{  1,2,\ldots,n\right\}  ^{2};\\k=\ell}}s_{k,i}%
s_{\ell,i}\underbrace{a_{k,\ell}}_{\substack{=a_{\ell,\ell}\\\text{(since
}k=\ell\text{)}}}\\
&  \ \ \ \ \ \ \ \ \ \ +\sum_{\substack{\left(  k,\ell\right)  \in\left\{
1,2,\ldots,n\right\}  ^{2};\\k>\ell}}\underbrace{s_{k,i}s_{\ell,i}}%
_{=s_{\ell,i}s_{k,i}}\underbrace{a_{k,\ell}}_{\substack{=-a_{\ell
,k}\\\text{(by Lemma \ref{lem.sol.altern.STAS.3} \textbf{(b)}}\\\text{(applied
to }\left(  k,\ell\right)  \text{ instead of }\left(  i,j\right)  \text{))}%
}}\\
&  \ \ \ \ \ \ \ \ \ \ \left(
\begin{array}
[c]{c}%
\text{since each }\left(  k,\ell\right)  \in\left\{  1,2,\ldots,n\right\}
^{2}\text{ satisfies exactly one of}\\
\text{the three assertions }\left(  k<\ell\right)  \text{, }\left(
k=\ell\right)  \text{ and }\left(  k>\ell\right)
\end{array}
\right) \\
&  =\underbrace{\sum_{\substack{\left(  k,\ell\right)  \in\left\{
1,2,\ldots,n\right\}  ^{2};\\k<\ell}}s_{k,i}s_{\ell,i}a_{k,\ell}%
}_{\substack{=\sum_{\substack{\left(  \ell,k\right)  \in\left\{
1,2,\ldots,n\right\}  ^{2};\\\ell<k}}s_{\ell,i}s_{k,i}a_{\ell,k}\\\text{(here,
we have renamed the}\\\text{summation index }\left(  k,\ell\right)  \text{ as
}\left(  \ell,k\right)  \text{)}}}+\sum_{\substack{\left(  k,\ell\right)
\in\left\{  1,2,\ldots,n\right\}  ^{2};\\k=\ell}}s_{k,i}s_{\ell,i}%
\underbrace{a_{\ell,\ell}}_{\substack{=0\\\text{(by Lemma
\ref{lem.sol.altern.STAS.3} \textbf{(a)}}\\\text{(applied to }\ell\text{
instead of }i\text{))}}}\\
&  \ \ \ \ \ \ \ \ \ \ +\underbrace{\sum_{\substack{\left(  k,\ell\right)
\in\left\{  1,2,\ldots,n\right\}  ^{2};\\k>\ell}}s_{\ell,i}s_{k,i}\left(
-a_{\ell,k}\right)  }_{=-\sum_{\substack{\left(  k,\ell\right)  \in\left\{
1,2,\ldots,n\right\}  ^{2};\\k>\ell}}s_{\ell,i}s_{k,i}a_{\ell,k}}\\
&  =\sum_{\substack{\left(  \ell,k\right)  \in\left\{  1,2,\ldots,n\right\}
^{2};\\\ell<k}}s_{\ell,i}s_{k,i}a_{\ell,k}+\underbrace{\sum_{\substack{\left(
k,\ell\right)  \in\left\{  1,2,\ldots,n\right\}  ^{2};\\k=\ell}}s_{k,i}%
s_{\ell,i}0}_{=0}+\left(  -\sum_{\substack{\left(  k,\ell\right)  \in\left\{
1,2,\ldots,n\right\}  ^{2};\\k>\ell}}s_{\ell,i}s_{k,i}a_{\ell,k}\right)
\end{align*}%
\begin{align*}
&  =\sum_{\substack{\left(  \ell,k\right)  \in\left\{  1,2,\ldots,n\right\}
^{2};\\\ell<k}}s_{\ell,i}s_{k,i}a_{\ell,k}-\underbrace{\sum_{\substack{\left(
k,\ell\right)  \in\left\{  1,2,\ldots,n\right\}  ^{2};\\k>\ell}}}%
_{\substack{=\sum_{\substack{\left(  \ell,k\right)  \in\left\{  1,2,\ldots
,n\right\}  ^{2};\\k>\ell}}=\sum_{\substack{\left(  \ell,k\right)  \in\left\{
1,2,\ldots,n\right\}  ^{2};\\\ell<k}}\\\text{(because for every }\left(
\ell,k\right)  \in\left\{  1,2,\ldots,n\right\}  ^{2}\text{,}\\\text{the
statement }\left(  k>\ell\right)  \text{ is equivalent to }\left(
\ell<k\right)  \text{)}}}s_{\ell,i}s_{k,i}a_{\ell,k}\\
&  =\sum_{\substack{\left(  \ell,k\right)  \in\left\{  1,2,\ldots,n\right\}
^{2};\\\ell<k}}s_{\ell,i}s_{k,i}a_{\ell,k}-\sum_{\substack{\left(
\ell,k\right)  \in\left\{  1,2,\ldots,n\right\}  ^{2};\\\ell<k}}s_{\ell
,i}s_{k,i}a_{\ell,k}=0.
\end{align*}

\end{vershort}

\begin{verlong}
Let $i\in\left\{  1,2,\ldots,m\right\}  $. Thus, $\left(  i,i\right)
\in\left\{  1,2,\ldots,m\right\}  ^{2}$. Hence, (\ref{sol.altern.STAS.4})
(applied to $\left(  i,i\right)  $ instead of $\left(  i,j\right)  $) yields%
\begin{align*}
b_{i,i}  &  =\sum_{\left(  k,\ell\right)  \in\left\{  1,2,\ldots,n\right\}
^{2}}s_{k,i}s_{\ell,i}a_{k,\ell}\\
&  =\sum_{\substack{\left(  k,\ell\right)  \in\left\{  1,2,\ldots,n\right\}
^{2};\\k<\ell}}s_{k,i}s_{\ell,i}a_{k,\ell}+\sum_{\substack{\left(
k,\ell\right)  \in\left\{  1,2,\ldots,n\right\}  ^{2};\\k=\ell}}s_{k,i}%
s_{\ell,i}\underbrace{a_{k,\ell}}_{\substack{=a_{\ell,\ell}\\\text{(since
}k=\ell\text{)}}}\\
&  \ \ \ \ \ \ \ \ \ \ +\sum_{\substack{\left(  k,\ell\right)  \in\left\{
1,2,\ldots,n\right\}  ^{2};\\k>\ell}}\underbrace{s_{k,i}s_{\ell,i}}%
_{=s_{\ell,i}s_{k,i}}\underbrace{a_{k,\ell}}_{\substack{=-a_{\ell
,k}\\\text{(by Lemma \ref{lem.sol.altern.STAS.3} \textbf{(b)}}\\\text{(applied
to }\left(  k,\ell\right)  \text{ instead of }\left(  i,j\right)  \text{))}%
}}\\
&  \ \ \ \ \ \ \ \ \ \ \left(
\begin{array}
[c]{c}%
\text{since each }\left(  k,\ell\right)  \in\left\{  1,2,\ldots,n\right\}
^{2}\text{ satisfies exactly one of}\\
\text{the three assertions }\left(  k<\ell\right)  \text{, }\left(
k=\ell\right)  \text{ and }\left(  k>\ell\right)
\end{array}
\right) \\
&  =\underbrace{\sum_{\substack{\left(  k,\ell\right)  \in\left\{
1,2,\ldots,n\right\}  ^{2};\\k<\ell}}s_{k,i}s_{\ell,i}a_{k,\ell}%
}_{\substack{=\sum_{\substack{\left(  \ell,k\right)  \in\left\{
1,2,\ldots,n\right\}  ^{2};\\\ell<k}}s_{\ell,i}s_{k,i}a_{\ell,k}\\\text{(here,
we have renamed the}\\\text{summation index }\left(  k,\ell\right)  \text{ as
}\left(  \ell,k\right)  \text{)}}}+\sum_{\substack{\left(  k,\ell\right)
\in\left\{  1,2,\ldots,n\right\}  ^{2};\\k=\ell}}s_{k,i}s_{\ell,i}%
\underbrace{a_{\ell,\ell}}_{\substack{=0\\\text{(by Lemma
\ref{lem.sol.altern.STAS.3} \textbf{(a)}}\\\text{(applied to }\ell\text{
instead of }i\text{))}}}\\
&  \ \ \ \ \ \ \ \ \ \ +\underbrace{\sum_{\substack{\left(  k,\ell\right)
\in\left\{  1,2,\ldots,n\right\}  ^{2};\\k>\ell}}s_{\ell,i}s_{k,i}\left(
-a_{\ell,k}\right)  }_{=-\sum_{\substack{\left(  k,\ell\right)  \in\left\{
1,2,\ldots,n\right\}  ^{2};\\k>\ell}}s_{\ell,i}s_{k,i}a_{\ell,k}}\\
&  =\sum_{\substack{\left(  \ell,k\right)  \in\left\{  1,2,\ldots,n\right\}
^{2};\\\ell<k}}s_{\ell,i}s_{k,i}a_{\ell,k}+\underbrace{\sum_{\substack{\left(
k,\ell\right)  \in\left\{  1,2,\ldots,n\right\}  ^{2};\\k=\ell}}s_{k,i}%
s_{\ell,i}0}_{=0}+\left(  -\sum_{\substack{\left(  k,\ell\right)  \in\left\{
1,2,\ldots,n\right\}  ^{2};\\k>\ell}}s_{\ell,i}s_{k,i}a_{\ell,k}\right) \\
&  =\sum_{\substack{\left(  \ell,k\right)  \in\left\{  1,2,\ldots,n\right\}
^{2};\\\ell<k}}s_{\ell,i}s_{k,i}a_{\ell,k}+\left(  -\sum_{\substack{\left(
k,\ell\right)  \in\left\{  1,2,\ldots,n\right\}  ^{2};\\k>\ell}}s_{\ell
,i}s_{k,i}a_{\ell,k}\right)
\end{align*}%
\begin{align*}
&  =\sum_{\substack{\left(  \ell,k\right)  \in\left\{  1,2,\ldots,n\right\}
^{2};\\\ell<k}}s_{\ell,i}s_{k,i}a_{\ell,k}-\underbrace{\sum_{\substack{\left(
k,\ell\right)  \in\left\{  1,2,\ldots,n\right\}  ^{2};\\k>\ell}}}%
_{\substack{=\sum_{k\in\left\{  1,2,\ldots,n\right\}  }\sum_{\substack{\ell
\in\left\{  1,2,\ldots,n\right\}  ;\\k>\ell}}\\=\sum_{\ell\in\left\{
1,2,\ldots,n\right\}  }\sum_{\substack{k\in\left\{  1,2,\ldots,n\right\}
;\\k>\ell}}\\=\sum_{\substack{\left(  \ell,k\right)  \in\left\{
1,2,\ldots,n\right\}  ^{2};\\k>\ell}}=\sum_{\substack{\left(  \ell,k\right)
\in\left\{  1,2,\ldots,n\right\}  ^{2};\\\ell<k}}\\\text{(because for every
}\left(  \ell,k\right)  \in\left\{  1,2,\ldots,n\right\}  ^{2}\text{,}%
\\\text{the statement }\left(  k>\ell\right)  \text{ is equivalent to }\left(
\ell<k\right)  \text{)}}}s_{\ell,i}s_{k,i}a_{\ell,k}\\
&  =\sum_{\substack{\left(  \ell,k\right)  \in\left\{  1,2,\ldots,n\right\}
^{2};\\\ell<k}}s_{\ell,i}s_{k,i}a_{\ell,k}-\sum_{\substack{\left(
\ell,k\right)  \in\left\{  1,2,\ldots,n\right\}  ^{2};\\\ell<k}}s_{\ell
,i}s_{k,i}a_{\ell,k}=0.
\end{align*}

\end{verlong}

Now, forget that we fixed $i$. We thus have shown that $\left(  b_{i,i}%
=0\text{ for all }i\in\left\{  1,2,\ldots,m\right\}  \right)  $.

Now, recall that $B=\left(  b_{i,j}\right)  _{1\leq i\leq m,\ 1\leq j\leq m}$.
Hence, the $m\times m$-matrix $B$ is alternating if and only if it satisfies
$B^{T}=-B$ and $\left(  b_{i,i}=0\text{ for all }i\in\left\{  1,2,\ldots
,m\right\}  \right)  $ (by the definition of \textquotedblleft
alternating\textquotedblright). Thus, the $m\times m$-matrix $B$ is
alternating (since it satisfies $B^{T}=-B$ and $\left(  b_{i,i}=0\text{ for
all }i\in\left\{  1,2,\ldots,m\right\}  \right)  $). Since $B=S^{T}AS$, this
rewrites as follows: The $m\times m$-matrix $S^{T}AS$ is alternating. This
solves Exercise \ref{exe.altern.STAS}.
\end{proof}

\subsection{Solution to Exercise \ref{exe.altern.det}}

Before we solve Exercise \ref{exe.altern.det}, let us show a combinatorial lemma:

\begin{lemma}
\label{lem.sol.altern.det.mod2}Let $n\in\mathbb{N}$. Let $\sigma\in S_{n}$ be
such that $\sigma^{-1}=\sigma$. Let $X=\left\{  i\in\left\{  1,2,\ldots
,n\right\}  \ \mid\ \sigma\left(  i\right)  =i\right\}  $. Then,%
\[
\left\vert X\right\vert \equiv n\operatorname{mod}2.
\]

\end{lemma}

\begin{vershort}
\begin{proof}
[Proof of Lemma \ref{lem.sol.altern.det.mod2}.]Define two further sets $Y$ and
$Z$ by
\[
Y=\left\{  i\in\left\{  1,2,\ldots,n\right\}  \ \mid\ \sigma\left(  i\right)
<i\right\}
\]
and%
\[
Z=\left\{  i\in\left\{  1,2,\ldots,n\right\}  \ \mid\ \sigma\left(  i\right)
>i\right\}  .
\]

Every $j\in Y$ satisfies $\sigma\left(  j\right)  \in Z$%
\ \ \ \ \footnote{\textit{Proof.} Let $j\in Y$. Then,%
\[
j\in Y=\left\{  i\in\left\{  1,2,\ldots,n\right\}  \ \mid\ \sigma\left(
i\right)  <i\right\}  .
\]
In other words, $j$ is an element $i$ of $\left\{  1,2,\ldots,n\right\}  $
satisfying $\sigma\left(  i\right)  <i$. In other words, $j$ is an element of
$\left\{  1,2,\ldots,n\right\}  $ and satisfies $\sigma\left(  j\right)  <j$.
\par
Clearly, $\sigma\left(  j\right)  $ is an element of $\left\{  1,2,\ldots
,n\right\}  $ and satisfies $\sigma\left(  \sigma\left(  j\right)  \right)
>\sigma\left(  j\right)  $ (since $\sigma\left(  j\right)
<j=\underbrace{\sigma^{-1}}_{=\sigma}\left(  \sigma\left(  j\right)  \right)
=\sigma\left(  \sigma\left(  j\right)  \right)  $). In other words,
$\sigma\left(  j\right)  $ is an element $i$ of $\left\{  1,2,\ldots
,n\right\}  $ satisfying $\sigma\left(  i\right)  >i$. In other words,
$\sigma\left(  j\right)  \in\left\{  i\in\left\{  1,2,\ldots,n\right\}
\ \mid\ \sigma\left(  i\right)  >i\right\}  $. This rewrites as $\sigma\left(
j\right)  \in Z$ (since $Z=\left\{  i\in\left\{  1,2,\ldots,n\right\}
\ \mid\ \sigma\left(  i\right)  >i\right\}  $). Qed.}. Thus, we can define a
map $\alpha:Y\rightarrow Z$ by
\[
\left(  \alpha\left(  j\right)  =\sigma\left(  j\right)
\ \ \ \ \ \ \ \ \ \ \text{for every }j\in Y\right)  .
\]
Consider this map $\alpha$. Similarly, define a map $\beta:Z\rightarrow Y$ by
\[
\left(  \beta\left(  j\right)  =\sigma\left(  j\right)
\ \ \ \ \ \ \ \ \ \ \text{for every }j\in Z\right)  .
\]

We have $\alpha\circ\beta=\operatorname*{id}$\ \ \ \ \footnote{\textit{Proof.}
Let $j\in Z$. The definition of $\beta$ yields $\beta\left(  j\right)
=\sigma\left(  j\right)  $. Comparing this with $\underbrace{\sigma^{-1}%
}_{=\sigma}\left(  j\right)  =\sigma\left(  j\right)  $, we obtain
$\beta\left(  j\right)  =\sigma^{-1}\left(  j\right)  $. Now,%
\begin{align*}
\left(  \alpha\circ\beta\right)  \left(  j\right)   &  =\alpha\left(
\beta\left(  j\right)  \right)  =\sigma\left(  \underbrace{\beta\left(
j\right)  }_{=\sigma^{-1}\left(  j\right)  }\right)
\ \ \ \ \ \ \ \ \ \ \left(  \text{by the definition of }\alpha\right) \\
&  =\sigma\left(  \sigma^{-1}\left(  j\right)  \right)  =j=\operatorname*{id}%
\left(  j\right)  .
\end{align*}
\par
Now, forget that we fixed $j$. We thus have shown that $\left(  \alpha
\circ\beta\right)  \left(  j\right)  =\operatorname*{id}\left(  j\right)  $
for each $j\in Z$. In other words, $\alpha\circ\beta=\operatorname*{id}$.
Qed.} and $\beta\circ\alpha=\operatorname*{id}$\ \ \ \ \footnote{for similar
reasons}. Thus, the two maps $\alpha$ and $\beta$ are mutually inverse. Hence,
the map $\alpha$ is invertible, i.e., is a bijection. Thus, there exists a
bijection $Y\rightarrow Z$ (namely, the map $\alpha$). Hence, $\left\vert
Y\right\vert =\left\vert Z\right\vert $.

Now,%
\[
\sum_{i\in\left\{  1,2,\ldots,n\right\}  }1=\underbrace{\left\vert \left\{
1,2,\ldots,n\right\}  \right\vert }_{=n}\cdot1=n\cdot1=n.
\]
Hence,%
\begin{align*}
n  &  =\sum_{i\in\left\{  1,2,\ldots,n\right\}  }1\\
&  =\underbrace{\sum_{\substack{i\in\left\{  1,2,\ldots,n\right\}
;\\\sigma\left(  i\right)  <i}}}_{\substack{=\sum_{i\in Y}\\\text{(since
}\left\{  i\in\left\{  1,2,\ldots,n\right\}  \ \mid\ \sigma\left(  i\right)
<i\right\}  =Y\text{)}}}1+\underbrace{\sum_{\substack{i\in\left\{
1,2,\ldots,n\right\}  ;\\\sigma\left(  i\right)  =i}}}_{\substack{=\sum_{i\in
X}\\\text{(since }\left\{  i\in\left\{  1,2,\ldots,n\right\}  \ \mid
\ \sigma\left(  i\right)  =i\right\}  =X\text{)}}}1\\
&  \ \ \ \ \ \ \ \ \ \ +\underbrace{\sum_{\substack{i\in\left\{
1,2,\ldots,n\right\}  ;\\\sigma\left(  i\right)  >i}}}_{\substack{=\sum_{i\in
Z}\\\text{(since }\left\{  i\in\left\{  1,2,\ldots,n\right\}  \ \mid
\ \sigma\left(  i\right)  >i\right\}  =Z\text{)}}}1\\
&  \ \ \ \ \ \ \ \ \ \ \left(
\begin{array}
[c]{c}%
\text{since each }i\in\left\{  1,2,\ldots,n\right\}  \text{ satisfies exactly
one of the three}\\
\text{statements }\left(  \sigma\left(  i\right)  <i\right)  \text{, }\left(
\sigma\left(  i\right)  =i\right)  \text{ and }\left(  \sigma\left(  i\right)
>i\right)
\end{array}
\right) \\
&  =\underbrace{\sum_{i\in Y}1}_{=\left\vert Y\right\vert \cdot1=\left\vert
Y\right\vert }+\underbrace{\sum_{i\in X}1}_{=\left\vert X\right\vert
\cdot1=\left\vert X\right\vert }+\underbrace{\sum_{i\in Z}1}_{=\left\vert
Z\right\vert \cdot1=\left\vert Z\right\vert }=\underbrace{\left\vert
Y\right\vert }_{=\left\vert Z\right\vert }+\left\vert X\right\vert +\left\vert
Z\right\vert \\
&  =\left\vert Z\right\vert +\left\vert X\right\vert +\left\vert Z\right\vert
=\left\vert X\right\vert +2\cdot\left\vert Z\right\vert \equiv\left\vert
X\right\vert \operatorname{mod}2.
\end{align*}
This proves Lemma \ref{lem.sol.altern.det.mod2}.
\end{proof}
\end{vershort}

\begin{verlong}
\begin{proof}
[Proof of Lemma \ref{lem.sol.altern.det.mod2}.]We have $X=\left\{
i\in\left\{  1,2,\ldots,n\right\}  \ \mid\ \sigma\left(  i\right)  =i\right\}
$. Thus,%
\begin{equation}
\sum_{i\in X}=\sum_{\substack{i\in\left\{  1,2,\ldots,n\right\}
;\\\sigma\left(  i\right)  =i}} \label{pf.lem.sol.altern.det.mod2.X}%
\end{equation}
(an equality between summation signs).

Define a further set $Y$ by $Y=\left\{  i\in\left\{  1,2,\ldots,n\right\}
\ \mid\ \sigma\left(  i\right)  <i\right\}  $. Thus,%
\begin{equation}
\sum_{i\in Y}=\sum_{\substack{i\in\left\{  1,2,\ldots,n\right\}
;\\\sigma\left(  i\right)  <i}} \label{pf.lem.sol.altern.det.mod2.Y}%
\end{equation}
(an equality between summation signs).

Define a further set $Z$ by $Z=\left\{  i\in\left\{  1,2,\ldots,n\right\}
\ \mid\ \sigma\left(  i\right)  >i\right\}  $. Thus,%
\begin{equation}
\sum_{i\in Z}=\sum_{\substack{i\in\left\{  1,2,\ldots,n\right\}
;\\\sigma\left(  i\right)  >i}} \label{pf.lem.sol.altern.det.mod2.Z}%
\end{equation}
(an equality between summation signs).

Every $j\in Y$ satisfies $\sigma\left(  j\right)  \in Z$%
\ \ \ \ \footnote{\textit{Proof.} Let $j\in Y$. Then,%
\[
j\in Y=\left\{  i\in\left\{  1,2,\ldots,n\right\}  \ \mid\ \sigma\left(
i\right)  <i\right\}  .
\]
In other words, $j$ is an element $i$ of $\left\{  1,2,\ldots,n\right\}  $
satisfying $\sigma\left(  i\right)  <i$. In other words, $j$ is an element of
$\left\{  1,2,\ldots,n\right\}  $ and satisfies $\sigma\left(  j\right)  <j$.
\par
Clearly, $\sigma\left(  j\right)  $ is an element of $\left\{  1,2,\ldots
,n\right\}  $ and satisfies $\sigma\left(  \sigma\left(  j\right)  \right)
>\sigma\left(  j\right)  $ (since $\sigma\left(  j\right)
<j=\underbrace{\sigma^{-1}}_{=\sigma}\left(  \sigma\left(  j\right)  \right)
=\sigma\left(  \sigma\left(  j\right)  \right)  $). In other words,
$\sigma\left(  j\right)  $ is an element $i$ of $\left\{  1,2,\ldots
,n\right\}  $ satisfying $\sigma\left(  i\right)  >i$. In other words,
$\sigma\left(  j\right)  \in\left\{  i\in\left\{  1,2,\ldots,n\right\}
\ \mid\ \sigma\left(  i\right)  >i\right\}  $. This rewrites as $\sigma\left(
j\right)  \in Z$ (since $Z=\left\{  i\in\left\{  1,2,\ldots,n\right\}
\ \mid\ \sigma\left(  i\right)  >i\right\}  $). Qed.}. Thus, we can define a
map $\alpha:Y\rightarrow Z$ by
\[
\left(  \alpha\left(  j\right)  =\sigma\left(  j\right)
\ \ \ \ \ \ \ \ \ \ \text{for every }j\in Y\right)  .
\]
Consider this map $\alpha$.

Every $j\in Z$ satisfies $\sigma\left(  j\right)  \in Y$%
\ \ \ \ \footnote{\textit{Proof.} Let $j\in Z$. Then,%
\[
j\in Z=\left\{  i\in\left\{  1,2,\ldots,n\right\}  \ \mid\ \sigma\left(
i\right)  >i\right\}  .
\]
In other words, $j$ is an element $i$ of $\left\{  1,2,\ldots,n\right\}  $
satisfying $\sigma\left(  i\right)  >i$. In other words, $j$ is an element of
$\left\{  1,2,\ldots,n\right\}  $ and satisfies $\sigma\left(  j\right)  >j$.
\par
Clearly, $\sigma\left(  j\right)  $ is an element of $\left\{  1,2,\ldots
,n\right\}  $ and satisfies $\sigma\left(  \sigma\left(  j\right)  \right)
<\sigma\left(  j\right)  $ (since $\sigma\left(  j\right)
>j=\underbrace{\sigma^{-1}}_{=\sigma}\left(  \sigma\left(  j\right)  \right)
=\sigma\left(  \sigma\left(  j\right)  \right)  $). In other words,
$\sigma\left(  j\right)  $ is an element $i$ of $\left\{  1,2,\ldots
,n\right\}  $ satisfying $\sigma\left(  i\right)  <i$. In other words,
$\sigma\left(  j\right)  \in\left\{  i\in\left\{  1,2,\ldots,n\right\}
\ \mid\ \sigma\left(  i\right)  <i\right\}  $. This rewrites as $\sigma\left(
j\right)  \in Y$ (since $Y=\left\{  i\in\left\{  1,2,\ldots,n\right\}
\ \mid\ \sigma\left(  i\right)  <i\right\}  $). Qed.}. Thus, we can define a
map $\beta:Z\rightarrow Y$ by
\[
\left(  \beta\left(  j\right)  =\sigma\left(  j\right)
\ \ \ \ \ \ \ \ \ \ \text{for every }j\in Z\right)  .
\]
Consider this map $\beta$.

We have $\alpha\circ\beta=\operatorname*{id}$\ \ \ \ \footnote{\textit{Proof.}
Let $j\in Z$. The definition of $\beta$ yields $\beta\left(  j\right)
=\sigma\left(  j\right)  $. Comparing this with $\underbrace{\sigma^{-1}%
}_{=\sigma}\left(  j\right)  =\sigma\left(  j\right)  $, we obtain
$\beta\left(  j\right)  =\sigma^{-1}\left(  j\right)  $. Now,%
\begin{align*}
\left(  \alpha\circ\beta\right)  \left(  j\right)   &  =\alpha\left(
\beta\left(  j\right)  \right)  =\sigma\left(  \underbrace{\beta\left(
j\right)  }_{=\sigma^{-1}\left(  j\right)  }\right)
\ \ \ \ \ \ \ \ \ \ \left(  \text{by the definition of }\alpha\right) \\
&  =\sigma\left(  \sigma^{-1}\left(  j\right)  \right)  =j=\operatorname*{id}%
\left(  j\right)  .
\end{align*}
\par
Now, forget that we fixed $j$. We thus have shown that $\left(  \alpha
\circ\beta\right)  \left(  j\right)  =\operatorname*{id}\left(  j\right)  $
for each $j\in Z$. In other words, $\alpha\circ\beta=\operatorname*{id}$.
Qed.} and $\beta\circ\alpha=\operatorname*{id}$%
\ \ \ \ \footnote{\textit{Proof.} Let $j\in Y$. The definition of $\alpha$
yields $\alpha\left(  j\right)  =\sigma\left(  j\right)  $. Comparing this
with $\underbrace{\sigma^{-1}}_{=\sigma}\left(  j\right)  =\sigma\left(
j\right)  $, we obtain $\alpha\left(  j\right)  =\sigma^{-1}\left(  j\right)
$. Now,%
\begin{align*}
\left(  \beta\circ\alpha\right)  \left(  j\right)   &  =\beta\left(
\alpha\left(  j\right)  \right)  =\sigma\left(  \underbrace{\alpha\left(
j\right)  }_{=\sigma^{-1}\left(  j\right)  }\right)
\ \ \ \ \ \ \ \ \ \ \left(  \text{by the definition of }\beta\right) \\
&  =\sigma\left(  \sigma^{-1}\left(  j\right)  \right)  =j=\operatorname*{id}%
\left(  j\right)  .
\end{align*}
\par
Now, forget that we fixed $j$. We thus have shown that $\left(  \beta
\circ\alpha\right)  \left(  j\right)  =\operatorname*{id}\left(  j\right)  $
for each $j\in Y$. In other words, $\beta\circ\alpha=\operatorname*{id}$.
Qed.}. Thus, the two maps $\alpha$ and $\beta$ are mutually inverse. Hence,
the map $\alpha$ is invertible, i.e., is a bijection. Thus, there exists a
bijection $Y\rightarrow Z$ (namely, the map $\alpha$). Hence, $\left\vert
Y\right\vert =\left\vert Z\right\vert $.

Now,%
\[
\sum_{i\in\left\{  1,2,\ldots,n\right\}  }1=\underbrace{\left\vert \left\{
1,2,\ldots,n\right\}  \right\vert }_{=n}\cdot1=n\cdot1=n.
\]
Hence,%
\begin{align*}
n  &  =\sum_{i\in\left\{  1,2,\ldots,n\right\}  }1=\underbrace{\sum
_{\substack{i\in\left\{  1,2,\ldots,n\right\}  ;\\\sigma\left(  i\right)
<i}}}_{\substack{=\sum_{i\in Y}\\\text{(by (\ref{pf.lem.sol.altern.det.mod2.Y}%
))}}}1+\underbrace{\sum_{\substack{i\in\left\{  1,2,\ldots,n\right\}
;\\\sigma\left(  i\right)  =i}}}_{\substack{=\sum_{i\in X}\\\text{(by
(\ref{pf.lem.sol.altern.det.mod2.X}))}}}1+\underbrace{\sum_{\substack{i\in
\left\{  1,2,\ldots,n\right\}  ;\\\sigma\left(  i\right)  >i}}}%
_{\substack{=\sum_{i\in Z}\\\text{(by (\ref{pf.lem.sol.altern.det.mod2.Z}))}%
}}1\\
&  \ \ \ \ \ \ \ \ \ \ \left(
\begin{array}
[c]{c}%
\text{since each }i\in\left\{  1,2,\ldots,n\right\}  \text{ satisfies exactly
one of the three}\\
\text{statements }\left(  \sigma\left(  i\right)  <i\right)  \text{, }\left(
\sigma\left(  i\right)  =i\right)  \text{ and }\left(  \sigma\left(  i\right)
>i\right)
\end{array}
\right) \\
&  =\underbrace{\sum_{i\in Y}1}_{=\left\vert Y\right\vert \cdot1=\left\vert
Y\right\vert }+\underbrace{\sum_{i\in X}1}_{=\left\vert X\right\vert
\cdot1=\left\vert X\right\vert }+\underbrace{\sum_{i\in Z}1}_{=\left\vert
Z\right\vert \cdot1=\left\vert Z\right\vert }=\underbrace{\left\vert
Y\right\vert }_{=\left\vert Z\right\vert }+\left\vert X\right\vert +\left\vert
Z\right\vert \\
&  =\left\vert Z\right\vert +\left\vert X\right\vert +\left\vert Z\right\vert
=\left\vert X\right\vert +2\cdot\left\vert Z\right\vert \equiv\left\vert
X\right\vert \operatorname{mod}2.
\end{align*}
This proves Lemma \ref{lem.sol.altern.det.mod2}.
\end{proof}
\end{verlong}

\begin{corollary}
\label{cor.sol.altern.det.sigma0}Let $n\in\mathbb{N}$ be odd. Let $\sigma\in
S_{n}$ be such that $\sigma^{-1}=\sigma$. Let $A=\left(  a_{i,j}\right)
_{1\leq i\leq n,\ 1\leq j\leq n}$ be an alternating $n\times n$-matrix. Then,%
\[
\prod_{i=1}^{n}a_{i,\sigma\left(  i\right)  }=0.
\]

\end{corollary}

\begin{proof}
[Proof of Corollary \ref{cor.sol.altern.det.sigma0}.]Let $X=\left\{
i\in\left\{  1,2,\ldots,n\right\}  \ \mid\ \sigma\left(  i\right)  =i\right\}
$. Then, Lemma \ref{lem.sol.altern.det.mod2} yields $\left\vert X\right\vert
\equiv n\equiv1\operatorname{mod}2$ (since $n$ is odd). If we had
$X=\varnothing$, then we would have $\left\vert \underbrace{X}_{=\varnothing
}\right\vert =\left\vert \varnothing\right\vert =0\not \equiv
1\operatorname{mod}2$; this would contradict $\left\vert X\right\vert
\equiv1\operatorname{mod}2$. Hence, we cannot have $X=\varnothing$. Thus, we
have $X\neq\varnothing$. In other words, the set $X$ is nonempty. Hence, there
exists some $x\in X$. Consider this $x$.

We have $x\in X=\left\{  i\in\left\{  1,2,\ldots,n\right\}  \ \mid
\ \sigma\left(  i\right)  =i\right\}  $. In other words, $x$ is an element $i$
of $\left\{  1,2,\ldots,n\right\}  $ satisfying $\sigma\left(  i\right)  =i$.
In other words, $x$ is an element of $\left\{  1,2,\ldots,n\right\}  $ and
satisfies $\sigma\left(  x\right)  =x$. Now, Lemma \ref{lem.sol.altern.STAS.3}
\textbf{(a)} (applied to $i=x$) yields $a_{x,x}=0$. Since $\sigma\left(
x\right)  =x$, we have $a_{x,\sigma\left(  x\right)  }=a_{x,x}=0$.

But $x\in\left\{  1,2,\ldots,n\right\}  $. Hence, $a_{x,\sigma\left(
x\right)  }$ is a factor of the product $\prod_{i=1}^{n}a_{i,\sigma\left(
i\right)  }$ (namely, the factor for $i=x$). This factor is $0$ (since
$a_{x,\sigma\left(  x\right)  }=0$). Hence, one factor of the product
$\prod_{i=1}^{n}a_{i,\sigma\left(  i\right)  }$ is $0$. Thus, the whole
product $\prod_{i=1}^{n}a_{i,\sigma\left(  i\right)  }$ must be $0$ (because
if one factor of a product is $0$, then the whole product must be $0$). In
other words, $\prod_{i=1}^{n}a_{i,\sigma\left(  i\right)  }=0$. This proves
Corollary \ref{cor.sol.altern.det.sigma0}.
\end{proof}

\begin{lemma}
\label{lem.sol.altern.det.sigma-1}Let $n\in\mathbb{N}$ be odd. Let $\sigma\in
S_{n}$. Let $A=\left(  a_{i,j}\right)  _{1\leq i\leq n,\ 1\leq j\leq n}$ be an
alternating $n\times n$-matrix. Then,%
\[
\prod_{i=1}^{n}a_{i,\sigma^{-1}\left(  i\right)  }=-\prod_{i=1}^{n}%
a_{i,\sigma\left(  i\right)  }.
\]

\end{lemma}

\begin{vershort}
\begin{proof}
[Proof of Lemma \ref{lem.sol.altern.det.sigma-1}.]We have $\sigma\in S_{n}$.
In other words, $\sigma$ is a permutation of the set $\left\{  1,2,\ldots
,n\right\}  $ (since $S_{n}$ is the set of all permutations of the set
$\left\{  1,2,\ldots,n\right\}  $). In other words, $\sigma$ is a bijection
$\left\{  1,2,\ldots,n\right\}  \rightarrow\left\{  1,2,\ldots,n\right\}  $.

Now,
\begin{align*}
\prod_{i=1}^{n}a_{i,\sigma^{-1}\left(  i\right)  }  &  =\prod_{i=1}%
^{n}\underbrace{a_{\sigma\left(  i\right)  ,\sigma^{-1}\left(  \sigma\left(
i\right)  \right)  }}_{\substack{=a_{\sigma\left(  i\right)  ,i}\\\text{(since
}\sigma^{-1}\left(  \sigma\left(  i\right)  \right)  =i\text{)}}}\\
&  \ \ \ \ \ \ \ \ \ \ \left(
\begin{array}
[c]{c}%
\text{here, we have substituted }\sigma\left(  i\right)  \text{ for }i\text{
in the product,}\\
\text{since }\sigma:\left\{  1,2,\ldots,n\right\}  \rightarrow\left\{
1,2,\ldots,n\right\}  \text{ is a bijection}%
\end{array}
\right) \\
&  =\prod_{i=1}^{n}\underbrace{a_{\sigma\left(  i\right)  ,i}}%
_{\substack{=-a_{i,\sigma\left(  i\right)  }\\\text{(by Lemma
\ref{lem.sol.altern.STAS.3} \textbf{(b)}}\\\text{(applied to }\left(
\sigma\left(  i\right)  ,i\right)  \text{ instead of }\left(  i,j\right)
\text{))}}}=\prod_{i=1}^{n}\left(  -a_{i,\sigma\left(  i\right)  }\right) \\
&  =\underbrace{\left(  -1\right)  ^{n}}_{\substack{=-1\\\text{(since }n\text{
is odd)}}}\prod_{i=1}^{n}a_{i,\sigma\left(  i\right)  }=-\prod_{i=1}%
^{n}a_{i,\sigma\left(  i\right)  }.
\end{align*}
This proves Lemma \ref{lem.sol.altern.det.sigma-1}.
\end{proof}
\end{vershort}

\begin{verlong}
\begin{proof}
[Proof of Lemma \ref{lem.sol.altern.det.sigma-1}.]We have $\sigma\in S_{n}$.
In other words, $\sigma$ is a permutation of the set $\left\{  1,2,\ldots
,n\right\}  $ (since $S_{n}$ is the set of all permutations of the set
$\left\{  1,2,\ldots,n\right\}  $). In other words, $\sigma$ is a bijection
$\left\{  1,2,\ldots,n\right\}  \rightarrow\left\{  1,2,\ldots,n\right\}  $.

We have%
\begin{equation}
a_{\sigma\left(  i\right)  ,i}=-a_{i,\sigma\left(  i\right)  }%
\ \ \ \ \ \ \ \ \ \ \text{for every }i\in\left\{  1,2,\ldots,n\right\}
\label{pf.lem.sol.altern.det.sigma-1.1}%
\end{equation}
\footnote{\textit{Proof of (\ref{pf.lem.sol.altern.det.sigma-1.1}):} Let
$i\in\left\{  1,2,\ldots,n\right\}  $. Then, $\sigma\left(  i\right)
\in\left\{  1,2,\ldots,n\right\}  $. Hence, $\left(  \sigma\left(  i\right)
,i\right)  \in\left\{  1,2,\ldots,n\right\}  ^{2}$ (since $\sigma\left(
i\right)  \in\left\{  1,2,\ldots,n\right\}  $ and $i\in\left\{  1,2,\ldots
,n\right\}  $). Thus, Lemma \ref{lem.sol.altern.STAS.3} \textbf{(b)} (applied
to $\left(  \sigma\left(  i\right)  ,i\right)  $ instead of $\left(
i,j\right)  $) yields $a_{\sigma\left(  i\right)  ,i}=-a_{i,\sigma\left(
i\right)  }$. This proves (\ref{pf.lem.sol.altern.det.sigma-1.1}).}.

Now,%
\begin{align*}
\underbrace{\prod_{i=1}^{n}}_{=\prod_{i\in\left\{  1,2,\ldots,n\right\}  }%
}a_{i,\sigma^{-1}\left(  i\right)  }  &  =\prod_{i\in\left\{  1,2,\ldots
,n\right\}  }a_{i,\sigma^{-1}\left(  i\right)  }=\prod_{i\in\left\{
1,2,\ldots,n\right\}  }\underbrace{a_{\sigma\left(  i\right)  ,\sigma
^{-1}\left(  \sigma\left(  i\right)  \right)  }}_{\substack{=a_{\sigma\left(
i\right)  ,i}\\\text{(since }\sigma^{-1}\left(  \sigma\left(  i\right)
\right)  =i\text{)}}}\\
&  \ \ \ \ \ \ \ \ \ \ \left(
\begin{array}
[c]{c}%
\text{here, we have substituted }\sigma\left(  i\right)  \text{ for }i\text{
in the product,}\\
\text{since }\sigma:\left\{  1,2,\ldots,n\right\}  \rightarrow\left\{
1,2,\ldots,n\right\}  \text{ is a bijection}%
\end{array}
\right) \\
&  =\prod_{i\in\left\{  1,2,\ldots,n\right\}  }\underbrace{a_{\sigma\left(
i\right)  ,i}}_{\substack{=-a_{i,\sigma\left(  i\right)  }\\\text{(by
(\ref{pf.lem.sol.altern.det.sigma-1.1}))}}}=\underbrace{\prod_{i\in\left\{
1,2,\ldots,n\right\}  }}_{=\prod_{i=1}^{n}}\underbrace{\left(  -a_{i,\sigma
\left(  i\right)  }\right)  }_{=\left(  -1\right)  a_{i,\sigma\left(
i\right)  }}\\
&  =\prod_{i=1}^{n}\left(  \left(  -1\right)  a_{i,\sigma\left(  i\right)
}\right)  =\underbrace{\left(  -1\right)  ^{n}}_{\substack{=-1\\\text{(since
}n\text{ is odd)}}}\prod_{i=1}^{n}a_{i,\sigma\left(  i\right)  }=\left(
-1\right)  \prod_{i=1}^{n}a_{i,\sigma\left(  i\right)  }\\
&  =-\prod_{i=1}^{n}a_{i,\sigma\left(  i\right)  }.
\end{align*}
This proves Lemma \ref{lem.sol.altern.det.sigma-1}.
\end{proof}
\end{verlong}

\begin{proof}
[Solution to Exercise \ref{exe.altern.det}.]\textbf{(a)} Assume that $A$ is antisymmetric.

Recall that the matrix $A$ is antisymmetric if and only if $A^{T}=-A$ (by the
definition of \textquotedblleft antisymmetric\textquotedblright). Hence,
$A^{T}=-A$ (since the matrix $A$ is antisymmetric).

Exercise \ref{exe.ps4.4} yields $\det\left(  A^{T}\right)  =\det A$. Hence,%
\begin{align*}
\det A  &  =\det\left(  \underbrace{A^{T}}_{=-A=\left(  -1\right)  A}\right)
=\det\left(  \left(  -1\right)  A\right)  =\underbrace{\left(  -1\right)
^{n}}_{\substack{=-1\\\text{(since }n\text{ is odd)}}}\det A\\
&  \ \ \ \ \ \ \ \ \ \ \left(  \text{by Proposition \ref{prop.det.scale}
(applied to }\lambda=-1\text{)}\right) \\
&  =\left(  -1\right)  \det A=-\det A.
\end{align*}
Adding $\det A$ to both sides of this equality, we obtain $\det A+\det A=0$.
Thus, $0=\det A+\det A=2\det A$. This solves Exercise \ref{exe.altern.det}
\textbf{(a)}.

\textbf{(b)} Assume that $A$ is alternating.

Write the $n\times n$-matrix $A$ in the form $A=\left(  a_{i,j}\right)
_{1\leq i\leq n,\ 1\leq j\leq n}$.

The set $S_{n}$ is finite. Hence, there exists a bijection $\beta
:S_{n}\rightarrow\left\{  1,2,\ldots,\left\vert S_{n}\right\vert \right\}  $.
Consider this bijection $\beta$. (Of course, $\left\vert S_{n}\right\vert
=n!$; but we will not use this.)

Recall that%
\begin{equation}
\left(  -1\right)  ^{\sigma^{-1}}=\left(  -1\right)  ^{\sigma}%
\ \ \ \ \ \ \ \ \ \ \text{for every }\sigma\in S_{n}.
\label{sol.altern.det.b.(-1)sigma}%
\end{equation}

\begin{vershort}
The map $S_{n}\rightarrow S_{n},\ \sigma\mapsto\sigma^{-1}$ is a bijection
(indeed, it is its own inverse).
\end{vershort}

\begin{verlong}
Define a map $\mathbf{j}:S_{n}\rightarrow S_{n}$ by $\left(  \mathbf{j}\left(
\sigma\right)  =\sigma^{-1}\text{ for every }\sigma\in S_{n}\right)  $. Then,
$\mathbf{j}$ is the map $S_{n}\rightarrow S_{n},\ \sigma\mapsto\sigma^{-1}$.

We have $\mathbf{j}\circ\mathbf{j}=\operatorname*{id}$%
\ \ \ \ \footnote{\textit{Proof.} Every $\sigma\in S_{n}$ satisfies%
\begin{align*}
\left(  \mathbf{j}\circ\mathbf{j}\right)  \left(  \sigma\right)   &
=\mathbf{j}\left(  \underbrace{\mathbf{j}\left(  \sigma\right)  }%
_{\substack{=\sigma^{-1}\\\text{(by the definition of }\mathbf{j}\text{)}%
}}\right)  =\mathbf{j}\left(  \sigma^{-1}\right)  =\left(  \sigma^{-1}\right)
^{-1}\ \ \ \ \ \ \ \ \ \ \left(  \text{by the definition of }\mathbf{j}\right)
\\
&  =\sigma=\operatorname*{id}\left(  \sigma\right)  .
\end{align*}
In other words, $\mathbf{j}\circ\mathbf{j}=\operatorname*{id}$. Qed.}. The two
maps $\mathbf{j}$ and $\mathbf{j}$ are mutually inverse (since $\mathbf{j}%
\circ\mathbf{j}=\operatorname*{id}$ and $\mathbf{j}\circ\mathbf{j}%
=\operatorname*{id}$). Hence, the map $\mathbf{j}$ is invertible. In other
words, the map $\mathbf{j}$ is a bijection. In other words, the map
$S_{n}\rightarrow S_{n},\ \sigma\mapsto\sigma^{-1}$ is a bijection (since the
map $\mathbf{j}$ is the map $S_{n}\rightarrow S_{n},\ \sigma\mapsto\sigma
^{-1}$).
\end{verlong}

Now,%
\begin{align}
&  \sum_{\substack{\sigma\in S_{n};\\\beta\left(  \sigma\right)  <\beta\left(
\sigma^{-1}\right)  }}\left(  -1\right)  ^{\sigma}\prod_{i=1}^{n}%
a_{i,\sigma\left(  i\right)  }\nonumber\\
&  =\underbrace{\sum_{\substack{\sigma\in S_{n};\\\beta\left(  \sigma
^{-1}\right)  <\beta\left(  \left(  \sigma^{-1}\right)  ^{-1}\right)  }%
}}_{\substack{=\sum_{\substack{\sigma\in S_{n};\\\beta\left(  \sigma
^{-1}\right)  <\beta\left(  \sigma\right)  }}\\\text{(since }\left(
\sigma^{-1}\right)  ^{-1}=\sigma\text{)}}}\underbrace{\left(  -1\right)
^{\sigma^{-1}}}_{\substack{=\left(  -1\right)  ^{\sigma}\\\text{(by
(\ref{sol.altern.det.b.(-1)sigma}))}}}\underbrace{\prod_{i=1}^{n}%
a_{i,\sigma^{-1}\left(  i\right)  }}_{\substack{=-\prod_{i=1}^{n}%
a_{i,\sigma\left(  i\right)  }\\\text{(by Lemma
\ref{lem.sol.altern.det.sigma-1})}}}\nonumber\\
&  \ \ \ \ \ \ \ \ \ \ \left(
\begin{array}
[c]{c}%
\text{here, we have substituted }\sigma^{-1}\text{ for }\sigma\text{ in the
sum,}\\
\text{since the map }S_{n}\rightarrow S_{n},\ \sigma\mapsto\sigma^{-1}\text{
is a bijection}%
\end{array}
\right) \nonumber\\
&  =\underbrace{\sum_{\substack{\sigma\in S_{n};\\\beta\left(  \sigma
^{-1}\right)  <\beta\left(  \sigma\right)  }}}_{\substack{=\sum
_{\substack{\sigma\in S_{n};\\\beta\left(  \sigma\right)  >\beta\left(
\sigma^{-1}\right)  }}\\\text{(since for every }\sigma\in S_{n}\text{,
the}\\\text{statement }\left(  \beta\left(  \sigma^{-1}\right)  <\beta\left(
\sigma\right)  \right)  \text{ is}\\\text{equivalent to }\left(  \beta\left(
\sigma\right)  >\beta\left(  \sigma^{-1}\right)  \right)  \text{)}}}\left(
-1\right)  ^{\sigma}\left(  -\prod_{i=1}^{n}a_{i,\sigma\left(  i\right)
}\right)  =\sum_{\substack{\sigma\in S_{n};\\\beta\left(  \sigma\right)
>\beta\left(  \sigma^{-1}\right)  }}\left(  -1\right)  ^{\sigma}\left(
-\prod_{i=1}^{n}a_{i,\sigma\left(  i\right)  }\right) \nonumber\\
&  =-\sum_{\substack{\sigma\in S_{n};\\\beta\left(  \sigma\right)
>\beta\left(  \sigma^{-1}\right)  }}\left(  -1\right)  ^{\sigma}\prod
_{i=1}^{n}a_{i,\sigma\left(  i\right)  }. \label{sol.altern.det.b.1}%
\end{align}

Also, every $\sigma\in S_{n}$ satisfying $\beta\left(  \sigma\right)
=\beta\left(  \sigma^{-1}\right)  $ must satisfy
\begin{equation}
\prod_{i=1}^{n}a_{i,\sigma\left(  i\right)  }=0 \label{sol.altern.det.b.2}%
\end{equation}
\footnote{\textit{Proof of (\ref{sol.altern.det.b.2}):} Let $\sigma\in S_{n}$
be such that $\beta\left(  \sigma\right)  =\beta\left(  \sigma^{-1}\right)  $.
\par
Recall that $\beta$ is a bijection. Hence, a map $\beta^{-1}$ is well-defined.
Now, $\beta^{-1}\left(  \underbrace{\beta\left(  \sigma\right)  }%
_{=\beta\left(  \sigma^{-1}\right)  }\right)  =\beta^{-1}\left(  \beta\left(
\sigma^{-1}\right)  \right)  =\sigma^{-1}$, so that $\sigma^{-1}=\beta
^{-1}\left(  \beta\left(  \sigma\right)  \right)  =\sigma$. Thus, Corollary
\ref{cor.sol.altern.det.sigma0} yields $\prod_{i=1}^{n}a_{i,\sigma\left(
i\right)  }=0$. This proves (\ref{sol.altern.det.b.2}).}.

Now, recall that $A=\left(  a_{i,j}\right)  _{1\leq i\leq n,\ 1\leq j\leq n}$.
Hence, the definition of $\det A$ yields%
\begin{align*}
&  \det A\\
&  =\sum_{\sigma\in S_{n}}\left(  -1\right)  ^{\sigma}\prod_{i=1}%
^{n}a_{i,\sigma\left(  i\right)  }\\
&  =\underbrace{\sum_{\substack{\sigma\in S_{n};\\\beta\left(  \sigma\right)
<\beta\left(  \sigma^{-1}\right)  }}\left(  -1\right)  ^{\sigma}\prod
_{i=1}^{n}a_{i,\sigma\left(  i\right)  }}_{\substack{=-\sum_{\substack{\sigma
\in S_{n};\\\beta\left(  \sigma\right)  >\beta\left(  \sigma^{-1}\right)
}}\left(  -1\right)  ^{\sigma}\prod_{i=1}^{n}a_{i,\sigma\left(  i\right)
}\\\text{(by (\ref{sol.altern.det.b.1}))}}}+\sum_{\substack{\sigma\in
S_{n};\\\beta\left(  \sigma\right)  =\beta\left(  \sigma^{-1}\right)
}}\left(  -1\right)  ^{\sigma}\underbrace{\prod_{i=1}^{n}a_{i,\sigma\left(
i\right)  }}_{\substack{=0\\\text{(by (\ref{sol.altern.det.b.2}))}}%
}+\sum_{\substack{\sigma\in S_{n};\\\beta\left(  \sigma\right)  >\beta\left(
\sigma^{-1}\right)  }}\left(  -1\right)  ^{\sigma}\prod_{i=1}^{n}%
a_{i,\sigma\left(  i\right)  }\\
&  \ \ \ \ \ \ \ \ \ \ \left(
\begin{array}
[c]{c}%
\text{since each }\sigma\in S_{n}\text{ satisfies exactly one of the three}\\
\text{statements }\left(  \beta\left(  \sigma\right)  <\beta\left(
\sigma^{-1}\right)  \right)  \text{, }\left(  \beta\left(  \sigma\right)
=\beta\left(  \sigma^{-1}\right)  \right)  \text{ and }\left(  \beta\left(
\sigma\right)  >\beta\left(  \sigma^{-1}\right)  \right)
\end{array}
\right) \\
&  =-\sum_{\substack{\sigma\in S_{n};\\\beta\left(  \sigma\right)
>\beta\left(  \sigma^{-1}\right)  }}\left(  -1\right)  ^{\sigma}\prod
_{i=1}^{n}a_{i,\sigma\left(  i\right)  }+\underbrace{\sum_{\substack{\sigma\in
S_{n};\\\beta\left(  \sigma\right)  =\beta\left(  \sigma^{-1}\right)
}}\left(  -1\right)  ^{\sigma}0}_{=0}+\sum_{\substack{\sigma\in S_{n}%
;\\\beta\left(  \sigma\right)  >\beta\left(  \sigma^{-1}\right)  }}\left(
-1\right)  ^{\sigma}\prod_{i=1}^{n}a_{i,\sigma\left(  i\right)  }\\
&  =-\sum_{\substack{\sigma\in S_{n};\\\beta\left(  \sigma\right)
>\beta\left(  \sigma^{-1}\right)  }}\left(  -1\right)  ^{\sigma}\prod
_{i=1}^{n}a_{i,\sigma\left(  i\right)  }+\sum_{\substack{\sigma\in
S_{n};\\\beta\left(  \sigma\right)  >\beta\left(  \sigma^{-1}\right)
}}\left(  -1\right)  ^{\sigma}\prod_{i=1}^{n}a_{i,\sigma\left(  i\right)  }=0.
\end{align*}
This solves Exercise \ref{exe.altern.det} \textbf{(b)}.
\end{proof}

A short solution of Exercise \ref{exe.altern.det} \textbf{(b)} using abstract
algebra appears in \cite[Theorem 12.2]{AndDosS}.

\subsection{Solution to Exercise \ref{exe.tridiag.isl}}

\begin{proof}
[Solution to Exercise \ref{exe.tridiag.isl}.]Define $n-1$ elements
$c_{1},c_{2},\ldots,c_{n-1}$ of $\mathbb{K}$ by%
\[
\left(  c_{i}=1\text{ for every }i\in\left\{  1,2,\ldots,n-1\right\}  \right)
.
\]
Define an $n\times n$-matrix $A$ as in Definition \ref{def.tridiag}. For every
two elements $x$ and $y$ of $\left\{  0,1,\ldots,n\right\}  $ satisfying
$x\leq y$, we define a $\left(  y-x\right)  \times\left(  y-x\right)  $-matrix
$A_{x,y}$ as in Proposition \ref{prop.tridiag.rec}.

Now, we claim that%
\begin{equation}
u_{i}=\det\left(  A_{0,i}\right)  \ \ \ \ \ \ \ \ \ \ \text{for every }%
i\in\left\{  0,1,\ldots,n\right\}  . \label{sol.tridiag.isl.u}%
\end{equation}

[\textit{Proof of (\ref{sol.tridiag.isl.u}):} We shall prove
(\ref{sol.tridiag.isl.u}) by strong induction over $i$. Thus, we fix some
$k\in\left\{  0,1,\ldots,n\right\}  $. We assume that (\ref{sol.tridiag.isl.u}%
) holds for all $i<k$. We now must show that (\ref{sol.tridiag.isl.u}) holds
for $i=k$. In other words, we must show that $u_{k}=\det\left(  A_{0,k}%
\right)  $.

This holds if $k\leq1$\ \ \ \ \footnote{\textit{Proof.} Proposition
\ref{prop.tridiag.rec} \textbf{(a)} (applied to $x=0$) yields $\det\left(
A_{0,0}\right)  =1$. Compared with $u_{0}=1$, this yields $u_{0}=\det\left(
A_{0,0}\right)  $.
\par
Proposition \ref{prop.tridiag.rec} \textbf{(b)} (applied to $x=0$) yields
$\det\left(  A_{0,1}\right)  =a_{1}$. Compared with $u_{1}=a_{1}$, this yields
$u_{1}=\det\left(  A_{0,1}\right)  $.
\par
Now, $u_{k}=\det\left(  A_{0,k}\right)  $ holds if $k=0$ (because we have
$u_{0}=\det\left(  A_{0,0}\right)  $), and also holds if $k=1$ (since
$u_{1}=\det\left(  A_{0,1}\right)  $). Therefore, $u_{k}=\det\left(
A_{0,k}\right)  $ holds if $k\leq1$. Qed.}. Hence, we can WLOG assume that we
don't have $k\leq1$. Assume this.

We have $k>1$ (since we don't have $k\leq1$). Thus, $k\geq2$. Hence, $k-1$ and
$k-2$ are nonnegative integers satisfying $k-1<k$ and $k-2<k$. Hence, we can
apply (\ref{sol.tridiag.isl.u}) to $i=k-1$ (since we have assumed that
(\ref{sol.tridiag.isl.u}) holds for all $i<k$). As a result, we obtain
$u_{k-1}=\det\left(  A_{0,k-1}\right)  $. Also, we can apply
(\ref{sol.tridiag.isl.u}) to $i=k-2$ (since we have assumed that
(\ref{sol.tridiag.isl.u}) holds for all $i<k$). As a result, we obtain
$u_{k-2}=\det\left(  A_{0,k-2}\right)  $.

We have $c_{k-1}=1$ (by the definition of $c_{k-1}$).

Now, $0\leq k-2$ (since $k\geq2$). Hence, Proposition \ref{prop.tridiag.rec}
\textbf{(c)} (applied to $x=0$ and $y=k$) yields%
\begin{align*}
\det\left(  A_{0,k}\right)   &  =a_{k}\det\left(  A_{0,k-1}\right)
-b_{k-1}\underbrace{c_{k-1}}_{=1}\det\left(  A_{0,k-2}\right) \\
&  =a_{k}\det\left(  A_{0,k-1}\right)  -b_{k-1}\det\left(  A_{0,k-2}\right)  .
\end{align*}
Comparing this with%
\begin{align*}
u_{k}  &  =a_{k}\underbrace{u_{k-1}}_{=\det\left(  A_{0,k-1}\right)  }%
-b_{k-1}\underbrace{u_{k-2}}_{=\det\left(  A_{0,k-2}\right)  }\\
&  \ \ \ \ \ \ \ \ \ \ \left(  \text{by the recursive definition of }\left(
u_{0},u_{1},\ldots,u_{n}\right)  \right) \\
&  =a_{k}\det\left(  A_{0,k-1}\right)  -b_{k-1}\det\left(  A_{0,k-2}\right)  ,
\end{align*}
we obtain $u_{k}=\det\left(  A_{0,k}\right)  $. In other words,
(\ref{sol.tridiag.isl.u}) holds for $i=k$. This completes the induction step.
Thus, (\ref{sol.tridiag.isl.u}) is proven.]

Next, we claim that%
\begin{equation}
v_{i}=\det\left(  A_{n-i,n}\right)  \ \ \ \ \ \ \ \ \ \ \text{for every }%
i\in\left\{  0,1,\ldots,n\right\}  . \label{sol.tridiag.isl.v}%
\end{equation}

[\textit{Proof of (\ref{sol.tridiag.isl.v}):} We shall prove
(\ref{sol.tridiag.isl.v}) by strong induction over $i$. Thus, we fix some
$k\in\left\{  0,1,\ldots,n\right\}  $. We assume that (\ref{sol.tridiag.isl.v}%
) holds for all $i<k$. We now must show that (\ref{sol.tridiag.isl.v}) holds
for $i=k$. In other words, we must show that $v_{k}=\det\left(  A_{n-k,n}%
\right)  $.

This holds if $k\leq1$\ \ \ \ \footnote{\textit{Proof.} Proposition
\ref{prop.tridiag.rec} \textbf{(a)} (applied to $x=n$) yields $\det\left(
A_{n,n}\right)  =1$. Thus, $\det\left(  A_{n-0,n}\right)  =\det\left(
A_{n,n}\right)  =1$. Compared with $v_{0}=1$, this yields $v_{0}=\det\left(
A_{n-0,n}\right)  $.
\par
Proposition \ref{prop.tridiag.rec} \textbf{(b)} (applied to $x=n-1$) yields
$\det\left(  A_{n-1,\left(  n-1\right)  +1}\right)  =a_{\left(  n-1\right)
+1}$. This rewrites as $\det\left(  A_{n-1,n}\right)  =a_{n}$ (since $\left(
n-1\right)  +1=n$). Compared with $v_{1}=a_{n}$, this yields $v_{1}%
=\det\left(  A_{n-1,n}\right)  $.
\par
Now, $v_{k}=\det\left(  A_{n-k,n}\right)  $ holds if $k=0$ (because we have
$v_{0}=\det\left(  A_{n-0,n}\right)  $), and also holds if $k=1$ (since
$v_{1}=\det\left(  A_{n-1,n}\right)  $). Therefore, $v_{k}=\det\left(
A_{n-k,n}\right)  $ holds if $k\leq1$. Qed.}. Hence, we can WLOG assume that
we don't have $k\leq1$. Assume this.

We have $k>1$ (since we don't have $k\leq1$). Thus, $k\geq2$. Hence, $k-1$ and
$k-2$ are nonnegative integers satisfying $k-1<k$ and $k-2<k$. Hence, we can
apply (\ref{sol.tridiag.isl.v}) to $i=k-1$ (since we have assumed that
(\ref{sol.tridiag.isl.v}) holds for all $i<k$). As a result, we obtain
$v_{k-1}=\det\left(  \underbrace{A_{n-\left(  k-1\right)  ,n}}_{=A_{n-k+1,n}%
}\right)  =\det\left(  A_{n-k+1,n}\right)  $. Also, we can apply
(\ref{sol.tridiag.isl.v}) to $i=k-2$ (since we have assumed that
(\ref{sol.tridiag.isl.v}) holds for all $i<k$). As a result, we obtain
$v_{k-2}=\det\left(  \underbrace{A_{n-\left(  k-2\right)  ,n}}_{=A_{n-k+2,n}%
}\right)  =\det\left(  A_{n-k+2,n}\right)  $.

We have $c_{n-k+1}=1$ (by the definition of $c_{n-k+1}$).

Now, $n-k\leq n-2$ (since $k\geq2$). Hence, Proposition \ref{prop.tridiag.rec}
\textbf{(d)} (applied to $x=n-k$ and $y=n$) yields%
\begin{align*}
\det\left(  A_{n-k,n}\right)   &  =a_{n-k+1}\det\left(  A_{n-k+1,n}\right)
-b_{n-k+1}\underbrace{c_{n-k+1}}_{=1}\det\left(  A_{n-k+2,n}\right) \\
&  =a_{n-k+1}\det\left(  A_{n-k+1,n}\right)  -b_{n-k+1}\det\left(
A_{n-k+2,n}\right)  .
\end{align*}
Comparing this with%
\begin{align*}
v_{k}  &  =a_{n-k+1}\underbrace{v_{k-1}}_{=\det\left(  A_{n-k+1,n}\right)
}-b_{n-k+1}\underbrace{v_{k-2}}_{=\det\left(  A_{n-k+2,n}\right)  }\\
&  \ \ \ \ \ \ \ \ \ \ \left(  \text{by the recursive definition of }\left(
v_{0},v_{1},\ldots,v_{n}\right)  \right) \\
&  =a_{n-k+1}\det\left(  A_{n-k+1,n}\right)  -b_{n-k+1}\det\left(
A_{n-k+2,n}\right)  ,
\end{align*}
we obtain $v_{k}=\det\left(  A_{n-k,n}\right)  $. In other words,
(\ref{sol.tridiag.isl.v}) holds for $i=k$. This completes the induction step.
Thus, (\ref{sol.tridiag.isl.v}) is proven.]

Now, we are almost done. In fact, applying (\ref{sol.tridiag.isl.u}) to $i=n$,
we obtain $u_{n}=\det\left(  A_{0,n}\right)  $. On the other hand, applying
(\ref{sol.tridiag.isl.v}) to $i=n$, we obtain $v_{n}=\det\left(
A_{n-n,n}\right)  =\det\left(  A_{0,n}\right)  $ (since $n-n=0$). Comparing
this with $u_{n}=\det\left(  A_{0,n}\right)  $, we obtain $u_{n}=v_{n}$.
Exercise \ref{exe.tridiag.isl} is solved.
\end{proof}

\subsection{Solution to Exercise \ref{exe.tridiag.cf}}

\begin{proof}
[Solution to Exercise \ref{exe.tridiag.cf}.]Assume that all denominators
appearing in Exercise \ref{exe.tridiag.cf} are invertible. For every
$k\in\left\{  1,2,\ldots,n\right\}  $, define an element $p_{k}$ of
$\mathbb{K}$ by%
\[
p_{k}=a_{k}-\dfrac{b_{k}c_{k}}{a_{k+1}-\dfrac{b_{k+1}c_{k+1}}{a_{k+2}%
-\dfrac{b_{k+2}c_{k+2}}{%
\begin{array}
[c]{ccc}%
a_{k+3}- &  & \\
& \ddots & \\
&  & -\dfrac{b_{n-2}c_{n-2}}{a_{n-1}-\dfrac{b_{n-1}c_{n-1}}{a_{n}}}%
\end{array}
}}}%
\]
\footnote{If $k=n$, then this should be interpreted as saying that
$p_{n}=a_{n}$.}. This definition of $p_{k}$ immediately gives a recursion:

\begin{itemize}
\item We have $p_{n}=a_{n}$.

\item For every $k\in\left\{  1,2,\ldots,n-1\right\}  $, we have%
\begin{equation}
p_{k}=a_{k}-\dfrac{b_{k}c_{k}}{p_{k+1}}. \label{sol.tridiag.cf.1}%
\end{equation}

\end{itemize}

Now, we shall show that
\begin{equation}
\det\left(  A_{n-k-1,n}\right)  =p_{n-k}\det\left(  A_{n-k,n}\right)
\ \ \ \ \ \ \ \ \ \ \text{for every }k\in\left\{  0,1,\ldots,n-1\right\}  .
\label{sol.tridiag.cf.main}%
\end{equation}

[\textit{Proof of (\ref{sol.tridiag.cf.main}):} We shall prove
(\ref{sol.tridiag.cf.main}) by induction over $k$:

\textit{Induction base:} It is easy to see that (\ref{sol.tridiag.cf.main})
holds for $k=0$\ \ \ \ \footnote{\textit{Proof.} Proposition
\ref{prop.tridiag.rec} \textbf{(a)} (applied to $x=n$) yields $\det\left(
A_{n,n}\right)  =1$. Hence, $\underbrace{p_{n-0}}_{=p_{n}=a_{n}}\det\left(
\underbrace{A_{n-0,n}}_{=A_{n,n}}\right)  =a_{n}$.
\par
Proposition \ref{prop.tridiag.rec} \textbf{(b)} (applied to $x=n-1$) yields
$\det\left(  A_{n-1,\left(  n-1\right)  +1}\right)  =a_{\left(  n-1\right)
+1}$. This rewrites as $\det\left(  A_{n-1,n}\right)  =a_{n}$ (since $\left(
n-1\right)  +1=n$). Thus, $\det\left(  \underbrace{A_{n-0-1,n}}_{=A_{n-1,n}%
}\right)  =\det\left(  A_{n-1,n}\right)  =a_{n}$. Comparing this with
$\underbrace{p_{n-0}}_{=p_{n}=a_{n}}\det\left(  \underbrace{A_{n-0,n}%
}_{=A_{n,n}}\right)  =a_{n}\underbrace{\det\left(  A_{n,n}\right)  }%
_{=1}=a_{n}$, we obtain $\det\left(  A_{n-0-1,n}\right)  =p_{n-0}\det\left(
A_{n-0,n}\right)  $. In other words, (\ref{sol.tridiag.cf.main}) holds for
$k=0$. Qed.}. This completes the induction base.

\textit{Induction step:} Let $K\in\left\{  0,1,\ldots,n-1\right\}  $ be
positive. Assume that (\ref{sol.tridiag.cf.main}) holds for $k=K-1$. We shall
show that (\ref{sol.tridiag.cf.main}) holds for $k=K$.

We have assumed that (\ref{sol.tridiag.cf.main}) holds for $k=K-1$. In other
words, we have $\det\left(  A_{n-\left(  K-1\right)  -1,n}\right)
=p_{n-\left(  K-1\right)  }\det\left(  A_{n-\left(  K-1\right)  ,n}\right)  $.
Since $n-\left(  K-1\right)  =n-K+1$, this rewrites as $\det\left(
A_{n-K+1-1,n}\right)  =p_{n-K+1}\det\left(  A_{n-K+1,n}\right)  $. Solving
this for $\det\left(  A_{n-K+1,n}\right)  $, we obtain%
\begin{align}
\det\left(  A_{n-K+1,n}\right)   &  =\dfrac{\det\left(  A_{n-K+1-1,n}\right)
}{p_{n-K+1}}=\dfrac{1}{p_{n-K+1}}\det\left(  \underbrace{A_{n-K+1-1,n}%
}_{=A_{n-K,n}}\right) \nonumber\\
&  =\dfrac{1}{p_{n-K+1}}\det\left(  A_{n-K,n}\right)  .
\label{sol.tridiag.cf.main.pf.2}%
\end{align}

We have $K\in\left\{  1,2,\ldots,n-1\right\}  $ (since $K$ is positive and
belongs to $\left\{  0,1,\ldots,n-1\right\}  $). Hence, applying
(\ref{sol.tridiag.cf.1}) to $k=n-K$, we obtain%
\begin{equation}
p_{n-K}=a_{n-K}-\dfrac{b_{n-K}c_{n-K}}{p_{n-K+1}}.
\label{sol.tridiag.cf.main.pf.4}%
\end{equation}

We have $K\leq n-1$ and thus $n-K-1\geq0$. Moreover, $K\geq1$ (since $K$ is
positive), thus $n-\underbrace{K}_{\geq1}-1\leq n-1-1=n-2$. Hence, Proposition
\ref{prop.tridiag.rec} \textbf{(d)} (applied to $x=n-K-1$ and $y=n$) shows
that%
\begin{align*}
&  \det\left(  A_{n-K-1,n}\right) \\
&  =\underbrace{a_{\left(  n-K-1\right)  +1}}_{=a_{n-K}}\det\left(
\underbrace{A_{\left(  n-K-1\right)  +1,n}}_{=A_{n-K,n}}\right)
-\underbrace{b_{\left(  n-K-1\right)  +1}}_{=b_{n-K}}\underbrace{c_{\left(
n-K-1\right)  +1}}_{=c_{n-K}}\det\left(  \underbrace{A_{\left(  n-K-1\right)
+2,n}}_{=A_{n-K+1,n}}\right) \\
&  =a_{n-K}\det\left(  A_{n-K,n}\right)  -b_{n-K}c_{n-K}\underbrace{\det
\left(  A_{n-K+1,n}\right)  }_{\substack{=\dfrac{1}{p_{n-K+1}}\det\left(
A_{n-K,n}\right)  \\\text{(by (\ref{sol.tridiag.cf.main.pf.2}))}}}\\
&  =a_{n-K}\det\left(  A_{n-K,n}\right)  -b_{n-K}c_{n-K}\cdot\dfrac
{1}{p_{n-K+1}}\det\left(  A_{n-K,n}\right) \\
&  =\underbrace{\left(  a_{n-K}-b_{n-K}c_{n-K}\cdot\dfrac{1}{p_{n-K+1}%
}\right)  }_{\substack{=a_{n-K}-\dfrac{b_{n-K}c_{n-K}}{p_{n-K+1}}%
=p_{n-K}\\\text{(by (\ref{sol.tridiag.cf.main.pf.4}))}}}\det\left(
A_{n-K,n}\right)  =p_{n-K}\det\left(  A_{n-K,n}\right)  .
\end{align*}
In other words, (\ref{sol.tridiag.cf.main}) holds for $k=K$. This completes
the induction step. Thus, (\ref{sol.tridiag.cf.main}) is proven by induction.]

Now, we can apply (\ref{sol.tridiag.cf.main}) to $k=n-1$. This gives us%
\[
\det\left(  A_{n-\left(  n-1\right)  -1,n}\right)  =\underbrace{p_{n-\left(
n-1\right)  }}_{=p_{1}}\det\left(  \underbrace{A_{n-\left(  n-1\right)  ,n}%
}_{=A_{1,n}}\right)  =p_{1}\det\left(  A_{1,n}\right)  .
\]
Since $n-\left(  n-1\right)  -1=0$, this rewrites as $\det\left(
A_{0,n}\right)  =p_{1}\det\left(  A_{1,n}\right)  $. Since $A_{0,n}=A$ (by
Proposition \ref{prop.tridiag.rec} \textbf{(e)}), this simplifies to $\det
A=p_{1}\det\left(  A_{1,n}\right)  $. Hence,%
\[
\dfrac{\det A}{\det\left(  A_{1,n}\right)  }=p_{1}=a_{1}-\dfrac{b_{1}c_{1}%
}{a_{2}-\dfrac{b_{2}c_{2}}{a_{3}-\dfrac{b_{3}c_{3}}{%
\begin{array}
[c]{ccc}%
a_{4}- &  & \\
& \ddots & \\
&  & -\dfrac{b_{n-2}c_{n-2}}{a_{n-1}-\dfrac{b_{n-1}c_{n-1}}{a_{n}}}%
\end{array}
}}}%
\]
(by the definition of $p_{1}$).
\end{proof}

\subsection{Solution to Exercise \ref{exe.tridiag.fib}}

\begin{proof}
[First solution to Exercise \ref{exe.tridiag.fib}.]We claim that%
\begin{equation}
\det\left(  A_{0,i}\right)  =f_{i+1}\ \ \ \ \ \ \ \ \ \ \text{for every }%
i\in\left\{  0,1,\ldots,n\right\}  . \label{sol.tridiag.fib.main}%
\end{equation}

[\textit{Proof of (\ref{sol.tridiag.fib.main}):} We shall prove
(\ref{sol.tridiag.fib.main}) by strong induction over $i$. Thus, we fix some
$k\in\left\{  0,1,\ldots,n\right\}  $. We assume that
(\ref{sol.tridiag.fib.main}) holds for all $i<k$. We now must show that
(\ref{sol.tridiag.fib.main}) holds for $i=k$. In other words, we must show
that $\det\left(  A_{0,k}\right)  =f_{k+1}$.

This holds if $k\leq1$\ \ \ \ \footnote{\textit{Proof.} Proposition
\ref{prop.tridiag.rec} \textbf{(a)} (applied to $x=0$) yields $\det\left(
A_{0,0}\right)  =1$. Compared with $f_{0+1}=f_{1}=1$, this yields $\det\left(
A_{0,0}\right)  =f_{0+1}$.
\par
Proposition \ref{prop.tridiag.rec} \textbf{(b)} (applied to $x=0$) yields
$\det\left(  A_{0,1}\right)  =a_{1}=1$. Compared with $f_{1+1}=f_{2}=1$, this
yields $\det\left(  A_{0,1}\right)  =f_{1+1}$.
\par
Now, $\det\left(  A_{0,k}\right)  =f_{k+1}$ holds if $k=0$ (because we have
$\det\left(  A_{0,0}\right)  =f_{0+1}$), and also holds if $k=1$ (since
$\det\left(  A_{0,1}\right)  =f_{1+1}$). Therefore, $\det\left(
A_{0,k}\right)  =f_{k+1}$ holds if $k\leq1$. Qed.}. Hence, we can WLOG assume
that we don't have $k\leq1$. Assume this.

We have $k>1$ (since we don't have $k\leq1$). Thus, $k\geq2$. Hence, $k-1$ and
$k-2$ are nonnegative integers satisfying $k-1<k$ and $k-2<k$. Hence, we can
apply (\ref{sol.tridiag.fib.main}) to $i=k-1$ (since we have assumed that
(\ref{sol.tridiag.fib.main}) holds for all $i<k$). As a result, we obtain
$\det\left(  A_{0,k-1}\right)  =f_{\left(  k-1\right)  +1}$. Also, we can
apply (\ref{sol.tridiag.fib.main}) to $i=k-2$ (since we have assumed that
(\ref{sol.tridiag.fib.main}) holds for all $i<k$). As a result, we obtain
$\det\left(  A_{0,k-2}\right)  =f_{\left(  k-2\right)  +1}$.

Now, $0\leq k-2$ (since $k\geq2$). Hence, Proposition \ref{prop.tridiag.rec}
\textbf{(c)} (applied to $x=0$ and $y=k$) yields%
\begin{align*}
\det\left(  A_{0,k}\right)   &  =\underbrace{a_{k}}_{\substack{=1\\\text{(by
the}\\\text{definition of }a_{k}\text{)}}}\underbrace{\det\left(
A_{0,k-1}\right)  }_{\substack{=f_{\left(  k-1\right)  +1}\\=f_{k}%
}}-\underbrace{b_{k-1}}_{\substack{=1\\\text{(by the}\\\text{definition of
}b_{k-1}\text{)}}}\underbrace{c_{k-1}}_{\substack{=-1\\\text{(by
the}\\\text{definition of }c_{k-1}\text{)}}}\underbrace{\det\left(
A_{0,k-2}\right)  }_{\substack{=f_{\left(  k-2\right)  +1}\\=f_{k-1}}}\\
&  =f_{k}-\left(  -1\right)  f_{k-1}=f_{k}+f_{k-1}.
\end{align*}
Comparing this with%
\[
f_{k+1}=f_{k}+f_{k-1}\ \ \ \ \ \ \ \ \ \ \left(  \text{by the recursive
definition of the Fibonacci numbers}\right)  ,
\]
we obtain $\det\left(  A_{0,k}\right)  =f_{k+1}$. In other words,
(\ref{sol.tridiag.fib.main}) holds for $i=k$. This completes the induction
step. Thus, (\ref{sol.tridiag.fib.main}) is proven.]

Now, applying (\ref{sol.tridiag.fib.main}) to $i=n$, we obtain $\det\left(
A_{0,n}\right)  =f_{n+1}$. Since $A_{0,n}=A$ (by Proposition
\ref{prop.tridiag.rec} \textbf{(e)}), we can rewrite this as $\det A=f_{n+1}$.
This solves Exercise \ref{exe.tridiag.fib}.
\end{proof}

\begin{proof}
[Second solution to Exercise \ref{exe.tridiag.fib} (sketched).]Here is a
different solution for Exercise \ref{exe.tridiag.fib}, which is far more
complicated than the previous one, but has the pedagogical advantage of
illuminating the connection between determinants and permutations, and the
combinatorics of the latter.

Exercise \ref{exe.ps2.2.3} (applied to $n+1$ instead of $n$) shows that
$f_{n+1}$ is the number of subsets $I$ of $\left\{  1,2,\ldots,n-1\right\}  $
such that no two elements of $I$ are consecutive. We shall refer to such
subsets $I$ as \textit{lacunar sets}.\footnote{We keep $n$ fixed, so a
\textquotedblleft lacunar set\textquotedblright\ will always be a subset of
$\left\{  1,2,\ldots,n-1\right\}  $.} Thus, $f_{n+1}$ is the number of all
lacunar sets.

(For example, if $n=5$, then the lacunar sets are $\varnothing$, $\left\{
1\right\}  $, $\left\{  2\right\}  $, $\left\{  3\right\}  $, $\left\{
4\right\}  $, $\left\{  1,3\right\}  $, $\left\{  1,4\right\}  $, and
$\left\{  2,4\right\}  $. Their number, unsurprisingly, is $8=f_{5+1}$.)

For any $\sigma\in S_{n}$, we define the following terminology:

\begin{itemize}
\item The \textit{excedances} of $\sigma$ are the elements $i\in\left\{
1,2,\ldots,n\right\}  $ satisfying $\sigma\left(  i\right)  >i$. For instance,
the permutation in $S_{7}$ written in one-line notation as $\left(
3,1,2,4,5,7,6\right)  $ has excedances $1$ and $6$.

\item We let $\operatorname*{Exced}\sigma$ denote the set of all excedances of
$\sigma$.

\item A permutation $\sigma\in S_{n}$ is said to be \textit{short-legged} if
every $i\in\left\{  1,2,\ldots,n\right\}  $ satisfies $\left\vert
\sigma\left(  i\right)  -i\right\vert \leq1$. For instance, the permutation in
$S_{7}$ written in one-line notation as $\left(  1,2,4,3,5,7,6\right)  $ is short-legged.
\end{itemize}

(Most of the terminology here is my own, tailored for this exercise; only the
notion of \textquotedblleft excedance\textquotedblright\ is standard. I chose
the name \textquotedblleft short-legged\textquotedblright\ because a
permutation $\sigma$ satisfying $\left\vert \sigma\left(  i\right)
-i\right\vert \leq1$ \textquotedblleft does not take $i$ very
far\textquotedblright.)

What does this all have to do with the exercise? Let us write our matrix $A$
in the form $A=\left(  a_{i,j}\right)  _{1\leq i\leq n,\ 1\leq j\leq n}$.
Thus, for every $\left(  i,j\right)  \in\left\{  1,2,\ldots,n\right\}  ^{2}$,
we have%
\begin{equation}
a_{i,j}=%
\begin{cases}
a_{i}, & \text{if }i=j;\\
b_{i}, & \text{if }i=j-1;\\
c_{j}, & \text{if }i=j+1;\\
0, & \text{otherwise}%
\end{cases}
=%
\begin{cases}
1, & \text{if }i=j;\\
1, & \text{if }i=j-1;\\
-1, & \text{if }i=j+1;\\
0, & \text{otherwise}%
\end{cases}
\label{sol.tridiag.fib.sol2.aij}%
\end{equation}
(since $a_{i}=1$, $b_{i}=1$ and $c_{j}=-1$). Notice that, as a consequence of
this equality, we have%
\begin{equation}
a_{i,j}=0\ \ \ \ \ \ \ \ \ \ \text{whenever }\left\vert i-j\right\vert >1.
\label{sol.tridiag.fib.sol2.locality}%
\end{equation}

Now, recall that $A=\left(  a_{i,j}\right)  _{1\leq i\leq n,\ 1\leq j\leq n}$,
so that $A^{T}=\left(  a_{j,i}\right)  _{1\leq i\leq n,\ 1\leq j\leq n}$. But
Exercise \ref{exe.ps4.4} yields $\det\left(  A^{T}\right)  =\det A$. Hence,%
\begin{equation}
\det A=\det\left(  A^{T}\right)  =\sum_{\sigma\in S_{n}}\left(  -1\right)
^{\sigma}\prod_{i=1}^{n}a_{\sigma\left(  i\right)  ,i}
\label{sol.tridiag.fib.sol2.detA=}%
\end{equation}
(by (\ref{eq.det.eq.2}), applied to $A^{T}$ and $a_{j,i}$ instead of $A$ and
$a_{i,j}$). This is an expression for $\det A$, but in order to get any
mileage out of it we need to simplify the terms $\prod_{i=1}^{n}%
a_{\sigma\left(  i\right)  ,i}$ for $\sigma\in S_{n}$. This turns out to
depend on whether the permutation $\sigma$ is short-legged or not:

\begin{itemize}
\item If a permutation $\sigma\in S_{n}$ is not short-legged, then%
\begin{equation}
\prod_{i=1}^{n}a_{\sigma\left(  i\right)  ,i}=0.
\label{sol.tridiag.fib.sol2.term1}%
\end{equation}

[\textit{Proof of (\ref{sol.tridiag.fib.sol2.term1}):} Let $\sigma\in S_{n}$
be not short-legged. Thus, there exists a $k\in\left\{  1,2,\ldots,n\right\}
$ satisfying $\left\vert \sigma\left(  k\right)  -k\right\vert >1$. The factor
$a_{\sigma\left(  k\right)  ,k}$ corresponding to this $k$ must be $0$
(because of (\ref{sol.tridiag.fib.sol2.locality})); this forces the whole
product $\prod_{i=1}^{n}a_{\sigma\left(  i\right)  ,i}$ to become $0$. Thus,
(\ref{sol.tridiag.fib.sol2.term1}) follows.]

\item If a permutation $\sigma\in S_{n}$ is short-legged, then%
\begin{equation}
\prod_{i=1}^{n}a_{\sigma\left(  i\right)  ,i}=\left(  -1\right)  ^{\left\vert
\operatorname*{Exced}\sigma\right\vert }. \label{sol.tridiag.fib.sol2.term2}%
\end{equation}

[\textit{Proof of (\ref{sol.tridiag.fib.sol2.term2}):} Let $\sigma\in S_{n}$
be short-legged. Thus, every $i\in\left\{  1,2,\ldots,n\right\}  $ satisfies
$\left\vert \sigma\left(  i\right)  -i\right\vert \leq1$. Consequently, for
every $i\in\left\{  1,2,\ldots,n\right\}  $, we have
\begin{align*}
a_{\sigma\left(  i\right)  ,i}  &  =%
\begin{cases}
1, & \text{if }\sigma\left(  i\right)  =i;\\
1, & \text{if }\sigma\left(  i\right)  =i-1;\\
-1, & \text{if }\sigma\left(  i\right)  =i+1;\\
0, & \text{otherwise}%
\end{cases}
\ \ \ \ \ \ \ \ \ \ \left(  \text{by (\ref{sol.tridiag.fib.sol2.aij})}\right)
\\
&  =%
\begin{cases}
1, & \text{if }\sigma\left(  i\right)  =i;\\
1, & \text{if }\sigma\left(  i\right)  =i-1;\\
-1, & \text{if }\sigma\left(  i\right)  =i+1
\end{cases}
\\
&  \ \ \ \ \ \ \ \ \ \ \left(
\begin{array}
[c]{c}%
\text{since the inequality }\left\vert \sigma\left(  i\right)  -i\right\vert
\leq1\text{ ensures that one of the}\\
\text{conditions }\sigma\left(  i\right)  =i\text{, }\sigma\left(  i\right)
=i-1\text{ and }\sigma\left(  i\right)  =i+1\text{ must hold}%
\end{array}
\right) \\
&  =%
\begin{cases}
1, & \text{if }\sigma\left(  i\right)  \leq i;\\
-1, & \text{if }\sigma\left(  i\right)  >i
\end{cases}
.
\end{align*}
Thus,%
\begin{align*}
\prod_{i=1}^{n}a_{\sigma\left(  i\right)  ,i}  &  =\prod_{i=1}^{n}%
\begin{cases}
1, & \text{if }\sigma\left(  i\right)  \leq i;\\
-1, & \text{if }\sigma\left(  i\right)  >i
\end{cases}
=\left(  -1\right)  ^{\left(  \text{the number of all }i\in\left\{
1,2,\ldots,n\right\}  \text{ satisfying }\sigma\left(  i\right)  >i\right)
}\\
&  =\left(  -1\right)  ^{\left\vert \left\{  i\in\left\{  1,2,\ldots
,n\right\}  \ \mid\ \sigma\left(  i\right)  >i\right\}  \right\vert }=\left(
-1\right)  ^{\left\vert \operatorname*{Exced}\sigma\right\vert };
\end{align*}
thus, (\ref{sol.tridiag.fib.sol2.term2}) is proven.]
\end{itemize}

Now, (\ref{sol.tridiag.fib.sol2.detA=}) becomes%
\begin{align}
\det A  &  =\sum_{\sigma\in S_{n}}\left(  -1\right)  ^{\sigma}\prod_{i=1}%
^{n}a_{\sigma\left(  i\right)  ,i}\nonumber\\
&  =\sum_{\substack{\sigma\in S_{n};\\\sigma\text{ is not short-legged}%
}}\left(  -1\right)  ^{\sigma}\underbrace{\prod_{i=1}^{n}a_{\sigma\left(
i\right)  ,i}}_{\substack{=0\\\text{(by (\ref{sol.tridiag.fib.sol2.term1}))}%
}}+\sum_{\substack{\sigma\in S_{n};\\\sigma\text{ is short-legged}}}\left(
-1\right)  ^{\sigma}\underbrace{\prod_{i=1}^{n}a_{\sigma\left(  i\right)  ,i}%
}_{\substack{=\left(  -1\right)  ^{\left\vert \operatorname*{Exced}%
\sigma\right\vert }\\\text{(by (\ref{sol.tridiag.fib.sol2.term2}))}%
}}\nonumber\\
&  =\underbrace{\sum_{\substack{\sigma\in S_{n};\\\sigma\text{ is not
short-legged}}}\left(  -1\right)  ^{\sigma}0}_{=0}+\sum_{\substack{\sigma\in
S_{n};\\\sigma\text{ is short-legged}}}\left(  -1\right)  ^{\sigma}\left(
-1\right)  ^{\left\vert \operatorname*{Exced}\sigma\right\vert }\nonumber\\
&  =\sum_{\substack{\sigma\in S_{n};\\\sigma\text{ is short-legged}}}\left(
-1\right)  ^{\sigma}\left(  -1\right)  ^{\left\vert \operatorname*{Exced}%
\sigma\right\vert }. \label{sol.tridiag.fib.sol2.det1}%
\end{align}
We still don't see how this connects to $f_{n+1}$, though. So let us relate
short-legged permutations to lacunar sets.

For any lacunar set $I$, we can define a permutation $\tau_{I}\in S_{n}$ by
the following rule:%
\[
\tau_{I}\left(  k\right)  =%
\begin{cases}
k+1, & \text{if }k\in I;\\
k-1, & \text{if }k-1\in I;\\
k, & \text{otherwise}%
\end{cases}
\ \ \ \ \ \ \ \ \ \ \text{for every }k\in\left\{  1,2,\ldots,n\right\}  .
\]
In other words, $\tau_{I}$ is the permutation of $\left\{  1,2,\ldots
,n\right\}  $ which interchanges every element $i$ of $I$ with its successor
$i+1$, while leaving all remaining elements unchanged. Make sure you
understand why $\tau_{I}$ is a well-defined map $\left\{  1,2,\ldots
,n\right\}  \rightarrow\left\{  1,2,\ldots,n\right\}  $\ \ \ \ \footnote{Here
are the two things you need to check:
\par
\begin{itemize}
\item The term $%
\begin{cases}
k+1, & \text{if }k\in I;\\
k-1, & \text{if }k-1\in I;\\
k, & \text{otherwise}%
\end{cases}
$ is unambiguous, because no $k$ satisfies both $k\in I$ and $k-1\in I$ at the
same time. (Here is where you use the lacunarity of $I$.)
\par
\item We have $%
\begin{cases}
k+1, & \text{if }k\in I;\\
k-1, & \text{if }k-1\in I;\\
k, & \text{otherwise}%
\end{cases}
\in\left\{  1,2,\ldots,n\right\}  $ for every $k\in\left\{  1,2,\ldots
,n\right\}  $. (Here you use $I\subseteq\left\{  1,2,\ldots,n-1\right\}  $. If
$I$ were only a subset of $\left\{  1,2,\ldots,n\right\}  $, then this would
fall outside of $\left\{  1,2,\ldots,n\right\}  $ for $k=n$.)
\end{itemize}
} and a permutation\footnote{Indeed, $\tau_{I}\circ\tau_{I}=\operatorname*{id}%
$, so that $\tau_{I}$ is its own inverse.}.

(For example, if $n=7$ and $I=\left\{  2,5\right\}  $, then $\tau_{I}=\left(
1,3,2,4,6,5,7\right)  $ in one-line notation.)

It is clear that the permutation $\tau_{I}$ is short-legged. Moreover, it
satisfies%
\begin{equation}
\operatorname*{Exced}\left(  \tau_{I}\right)  =I;
\label{sol.tridiag.fib.sol2.tauI.Exced}%
\end{equation}
as a consequence, it is possible to reconstruct $I$ from $\tau_{I}$. Thus, the
permutations $\tau_{I}$ for distinct $I$ are distinct.

What is the sign $\left(  -1\right)  ^{\tau_{I}}$ ? It is easy to see (from
the construction of $\tau_{I}$) that the only inversions of $\tau_{I}$ are the
pairs $\left(  i,i+1\right)  $ for $i\in I$ (essentially, the short-leggedness
of $\tau_{I}$ prevents $\tau_{I}$ from changing the order of two non-adjacent
integers). Thus, the number of these inversions is $\left\vert I\right\vert $.
Thus, $\ell\left(  \tau_{I}\right)  =\left\vert I\right\vert $. Hence,
\begin{equation}
\left(  -1\right)  ^{\tau_{I}}=\left(  -1\right)  ^{\ell\left(  \tau
_{I}\right)  }=\left(  -1\right)  ^{\left\vert I\right\vert }.
\label{sol.tridiag.fib.sol2.tauI.sign}%
\end{equation}

So there are at least some short-legged permutations that we understand well:
the $\tau_{I}$ for lacunar sets $I$. Are there others?

It turns out that there aren't. Indeed,
\begin{equation}
\text{every short-legged }\sigma\in S_{n}\text{ has the form }\tau_{I}\text{
for some lacunar set }I. \label{sol.tridiag.fib.sol2.tauI.surj}%
\end{equation}
Before we can prove this, we shall prove two auxiliary observations:

\textit{Observation 1:} Let $\sigma\in S_{n}$ be short-legged. If
$i\in\left\{  1,2,\ldots,n\right\}  $ be such that $\sigma\left(  i\right)
=i+1$, then $\sigma\left(  i+1\right)  =i$.

\textit{Observation 2:} Let $\sigma\in S_{n}$ be short-legged. If
$i\in\left\{  1,2,\ldots,n\right\}  $ be such that $\sigma\left(  i\right)
=i-1$, then $\sigma\left(  i-1\right)  =i$.

[\textit{Proof of Observation 1.} Assume the contrary. Thus, there exists some
$i\in\left\{  1,2,\ldots,n\right\}  $ such that $\sigma\left(  i\right)  =i+1$
but $\sigma\left(  i+1\right)  \neq i$. We call such $i$'s \textit{evil}. By
our assumption, there exists at least one evil $i$. Consider the highest evil
$i$. Thus, $i+1$ is not evil.

Since $i$ is evil, we have $\sigma\left(  i\right)  =i+1$ but $\sigma\left(
i+1\right)  \neq i$. In particular, $i+1=\sigma\left(  i\right)  \in\left\{
1,2,\ldots,n\right\}  $, so that $\sigma\left(  i+1\right)  $ is well-defined.
Since $\sigma$ is short-legged, we have $\left\vert \sigma\left(  i+1\right)
-\left(  i+1\right)  \right\vert \leq1$. Hence, $\sigma\left(  i+1\right)  $
is either $i$ or $i+1$ or $i+2$. But $\sigma\left(  i+1\right)  $ cannot be
$i$ (since $\sigma\left(  i+1\right)  \neq i$) and cannot be $i+1$ either
(since this would cause $\sigma\left(  i+1\right)  =i+1=\sigma\left(
i\right)  $, which would contradict the injectivity of $\sigma$). Hence,
$\sigma\left(  i+1\right)  $ must be $i+2$. In other words, $\sigma\left(
i+1\right)  =i+2$. Moreover, the injectivity of $\sigma$ shows that
$\sigma\left(  i+2\right)  \neq\sigma\left(  i\right)  =i+1$, so that $i+1$ is
evil. But this contradicts the fact that $i+1$ is not evil. Thus, Observation
1 is proven.]

[\textit{Proof of Observation 2.} Analogous to Observation 1 (this time, take
the lowest evil $i$), and left to the reader.]

[\textit{Proof of (\ref{sol.tridiag.fib.sol2.tauI.surj}):} Let $\sigma\in
S_{n}$ be short-legged. We must show that $\sigma=\tau_{I}$ for some lacunar
set $I$.

We set $I=\operatorname*{Exced}\sigma$. (This is the only choice we can make
to have any hope for $\sigma=\tau_{I}$ to be true; indeed,
(\ref{sol.tridiag.fib.sol2.tauI.Exced}) ensures that if $\sigma=\tau_{I}$,
then $\operatorname*{Exced}\sigma=I$.)

We notice that%
\begin{equation}
\sigma\left(  i\right)  =i+1\ \ \ \ \ \ \ \ \ \ \text{for every }i\in I
\label{sol.tridiag.fib.sol2.tauI.surj.pf.1}%
\end{equation}
\footnote{\textit{Proof of (\ref{sol.tridiag.fib.sol2.tauI.surj.pf.1}):} Let
$i\in I$. Thus, $i\in I=\operatorname*{Exced}\sigma$. In other words, $i$ is
an excedance of $\sigma$. Hence, $\sigma\left(  i\right)  >i$. Since
$\left\vert \sigma\left(  i\right)  -i\right\vert \leq1$ (because $\sigma$ is
short-legged), this means that $\sigma\left(  i\right)  =i+1$, qed.}. Thus,
$n$ cannot belong to $I$ (since this would entail $\sigma\left(  n\right)
=n+1$, but $n+1\notin\left\{  1,2,\ldots,n\right\}  $). Hence, $I\subseteq
\left\{  1,2,\ldots,n-1\right\}  $.

Let us first show that $I$ is a lacunar set. Indeed, assume (for the sake of
contradiction) that this is not so. Then, there exists some $i\in I$ such that
$i+1\in I$. Consider such an $i$. We have $\sigma\left(  i\right)  =i+1$ (by
(\ref{sol.tridiag.fib.sol2.tauI.surj.pf.1})), and thus $\sigma\left(
i+1\right)  =i$ (by Observation 1). But $\sigma\left(  i+1\right)  =i+2$ (by
(\ref{sol.tridiag.fib.sol2.tauI.surj.pf.1}), applied to $i+1$ instead of $i$).
Hence, $i+2=\sigma\left(  i+1\right)  =i$, which is absurd. Hence, we have
found a contradiction. This finishes our proof that $I$ is a lacunar set.

We still need to show that we actually have $\sigma=\tau_{I}$. In other words,
we need to show that $\sigma\left(  k\right)  =\tau_{I}\left(  k\right)  $ for
every $k\in\left\{  1,2,\ldots,n\right\}  $.

So let us fix $k\in\left\{  1,2,\ldots,n\right\}  $, and let us show that
$\sigma\left(  k\right)  =\tau_{I}\left(  k\right)  $. We are in one of the
following three cases:

\textit{Case 1:} We have $k\in I$.

\textit{Case 2:} We have $k-1\in I$.

\textit{Case 3:} Neither $k\in I$ nor $k-1\in I$.

\begin{itemize}
\item Let us first consider Case 1. In this case, $k\in I$. Hence,
(\ref{sol.tridiag.fib.sol2.tauI.surj.pf.1}) (applied to $i=k$) yields
$\sigma\left(  k\right)  =k+1$. On the other hand, the definition of $\tau
_{I}$ shows that $\tau_{I}\left(  k\right)  =k+1$ as well. Thus,
$\sigma\left(  k\right)  =k+1=\tau_{I}\left(  k\right)  $. Hence,
$\sigma\left(  k\right)  =\tau_{I}\left(  k\right)  $ is proven in Case 1.

\item Let us now consider Case 2. In this case, $k-1\in I$. Hence,
(\ref{sol.tridiag.fib.sol2.tauI.surj.pf.1}) (applied to $i=k-1$) yields
$\sigma\left(  k-1\right)  =k$. Since $\sigma$ is injective, we have
$\sigma\left(  k\right)  \neq\sigma\left(  k-1\right)  =k$. Also, $I$ is
lacunar, so that $k-1\in I$ entails $k\notin I$; thus, $k\notin
I=\operatorname*{Exced}\sigma$, so that $k$ is not an excedance of $\sigma$.
In other words, $\sigma\left(  k\right)  \leq k$. Combined with $\sigma\left(
k\right)  \neq k$, this yields $\sigma\left(  k\right)  <k$. Since $\left\vert
\sigma\left(  k\right)  -k\right\vert \leq1$ (because $\sigma$ is
short-legged), this shows that $\sigma\left(  k\right)  =k-1$. On the other
hand, $\tau_{I}\left(  k\right)  =k-1$ by the definition of $\tau_{I}$. Thus,
$\sigma\left(  k\right)  =k-1=\tau_{I}\left(  k\right)  $. Hence,
$\sigma\left(  k\right)  =\tau_{I}\left(  k\right)  $ is proven in Case 2.

\item Let us finally consider Case 3. In this case, neither $k\in I$ nor
$k-1\in I$. Let us now show that $\sigma\left(  k\right)  =k$. Indeed, assume
the contrary. Thus, $\sigma\left(  k\right)  \neq k$. As in Case 2, we can use
this (and $k\notin I$) to show that $\sigma\left(  k\right)  =k-1$.
Observation 2 thus shows that $\sigma\left(  k-1\right)  =k>k-1$, so that
$k-1$ is an excedance of $\sigma$. In other words, $k-1\in
\operatorname*{Exced}\sigma=I$. This contradicts the assumption that we do not
have $k-1\in I$. This contradiction concludes our proof of $\sigma\left(
k\right)  =k$. On the other hand, $\tau_{I}\left(  k\right)  =k$ by the
definition of $\tau_{I}$. Thus, $\sigma\left(  k\right)  =k=\tau_{I}\left(
k\right)  $. Hence, $\sigma\left(  k\right)  =\tau_{I}\left(  k\right)  $ is
proven in Case 3.
\end{itemize}

We now have shown that $\sigma\left(  k\right)  =\tau_{I}\left(  k\right)  $
in all possible cases. Thus, $\sigma=\tau_{I}$. Since $I$ is a lacunar set,
this proves (\ref{sol.tridiag.fib.sol2.tauI.surj}).]

All we now need to do is combine our results. For every lacunar set $I$, we
have defined a short-legged permutation $\tau_{I}\in S_{n}$. Conversely, we
know (from (\ref{sol.tridiag.fib.sol2.tauI.surj})) that every short-legged
$\sigma\in S_{n}$ has the form $\tau_{I}$ for some lacunar set $I$; we also
know that this $I$ is uniquely determined by the $\sigma$ (since the
permutations $\tau_{I}$ for distinct $I$ are distinct). Thus, we have a
bijection between the lacunar sets and the short-legged permutations in
$S_{n}$; the bijection sends every $I$ to $\tau_{I}$. Consequently,%
\begin{align*}
&  \sum_{\substack{\sigma\in S_{n};\\\sigma\text{ is short-legged}}}\left(
-1\right)  ^{\sigma}\left(  -1\right)  ^{\left\vert \operatorname*{Exced}%
\sigma\right\vert }\\
&  =\sum_{I\text{ is a lacunar set}}\underbrace{\left(  -1\right)  ^{\tau_{I}%
}}_{\substack{=\left(  -1\right)  ^{\left\vert I\right\vert }\\\text{(by
(\ref{sol.tridiag.fib.sol2.tauI.sign}))}}}\underbrace{\left(  -1\right)
^{\left\vert \operatorname*{Exced}\left(  \tau_{I}\right)  \right\vert }%
}_{\substack{=\left(  -1\right)  ^{\left\vert I\right\vert }\\\text{(by
(\ref{sol.tridiag.fib.sol2.tauI.Exced}))}}}\\
&  =\sum_{I\text{ is a lacunar set}}\underbrace{\left(  -1\right)
^{\left\vert I\right\vert }\left(  -1\right)  ^{\left\vert I\right\vert }%
}_{=\left(  \left(  -1\right)  ^{\left\vert I\right\vert }\right)  ^{2}%
=1}=\sum_{I\text{ is a lacunar set}}1\\
&  =\left(  \text{the number of all lacunar sets}\right)  =f_{n+1}.
\end{align*}
Combining this with (\ref{sol.tridiag.fib.sol2.det1}), we conclude that $\det
A=f_{n+1}$.
\end{proof}

\subsection{Solution to Exercise \ref{exe.block2x2.mult}}

\begin{vershort}
\begin{proof}
[Solution to Exercise \ref{exe.block2x2.mult}.]This is another case where the
solution is really clear with the appropriate amount of waving hands and
pointing fingers, but on paper becomes nearly impossible to convey. I shall
therefore resort to formalism and computation.

Write the matrices $A$, $B$, $C$, $D$, $A^{\prime}$, $B^{\prime}$, $C^{\prime
}$ and $D^{\prime}$ in the forms%
\begin{align*}
A  &  =\left(  a_{i,j}\right)  _{1\leq i\leq n,\ 1\leq j\leq m}%
,\ \ \ \ \ \ \ \ \ \ B=\left(  b_{i,j}\right)  _{1\leq i\leq n,\ 1\leq j\leq
m^{\prime}},\\
C  &  =\left(  c_{i,j}\right)  _{1\leq i\leq n^{\prime},\ 1\leq j\leq
m},\ \ \ \ \ \ \ \ \ \ D=\left(  d_{i,j}\right)  _{1\leq i\leq n^{\prime
},\ 1\leq j\leq m^{\prime}},\\
A^{\prime}  &  =\left(  a_{i,j}^{\prime}\right)  _{1\leq i\leq m,\ 1\leq
j\leq\ell},\ \ \ \ \ \ \ \ \ \ B^{\prime}=\left(  b_{i,j}^{\prime}\right)
_{1\leq i\leq m,\ 1\leq j\leq\ell^{\prime}},\\
C^{\prime}  &  =\left(  c_{i,j}^{\prime}\right)  _{1\leq i\leq m^{\prime
},\ 1\leq j\leq\ell},\ \ \ \ \ \ \ \ \ \ D^{\prime}=\left(  d_{i,j}^{\prime
}\right)  _{1\leq i\leq m^{\prime},\ 1\leq j\leq\ell^{\prime}}.
\end{align*}
The definition of the $\left(  n+n^{\prime}\right)  \times\left(  m+m^{\prime
}\right)  $-matrix $\left(
\begin{array}
[c]{cc}%
A & B\\
C & D
\end{array}
\right)  $ shows that\footnote{Here and in the following, we use the symbol
\textquotedblleft\&\textquotedblright\ as shorthand for the word
\textquotedblleft and\textquotedblright.}%
\begin{align*}
\left(
\begin{array}
[c]{cc}%
A & B\\
C & D
\end{array}
\right)   &  =\left(
\begin{array}
[c]{cccccccc}%
a_{1,1} & a_{1,2} & \cdots & a_{1,m} & b_{1,1} & b_{1,2} & \cdots &
b_{1,m^{\prime}}\\
a_{2,1} & a_{2,2} & \cdots & a_{2,m} & b_{2,1} & b_{2,2} & \cdots &
b_{2,m^{\prime}}\\
\vdots & \vdots & \ddots & \vdots & \vdots & \vdots & \ddots & \vdots\\
a_{n,1} & a_{n,2} & \cdots & a_{n,m} & b_{n,1} & b_{n,2} & \cdots &
b_{n,m^{\prime}}\\
c_{1,1} & c_{1,2} & \cdots & c_{1,m} & d_{1,1} & d_{1,2} & \cdots &
d_{1,m^{\prime}}\\
c_{2,1} & c_{2,2} & \cdots & c_{2,m} & d_{2,1} & d_{2,2} & \cdots &
d_{2,m^{\prime}}\\
\vdots & \vdots & \ddots & \vdots & \vdots & \vdots & \ddots & \vdots\\
c_{n^{\prime},1} & c_{n^{\prime},2} & \cdots & c_{n^{\prime},m} &
d_{n^{\prime},1} & d_{n^{\prime},2} & \cdots & d_{n^{\prime},m^{\prime}}%
\end{array}
\right) \\
&  =\left(
\begin{cases}
a_{i,j}, & \text{if }i\leq n\text{ \& }j\leq m;\\
b_{i,j-m}, & \text{if }i\leq n\text{ \& }j>m;\\
c_{i-n,j}, & \text{if }i>n\text{ \& }j\leq m;\\
d_{i-n,j-m}, & \text{if }i>n\text{ \& }j>m
\end{cases}
\right)  _{1\leq i\leq n+n^{\prime},\ 1\leq j\leq m+m^{\prime}}.
\end{align*}
Similarly,%
\begin{align*}
\left(
\begin{array}
[c]{cc}%
A^{\prime} & B^{\prime}\\
C^{\prime} & D^{\prime}%
\end{array}
\right)   &  =\left(
\begin{array}
[c]{cccccccc}%
a_{1,1}^{\prime} & a_{1,2}^{\prime} & \cdots & a_{1,\ell}^{\prime} &
b_{1,1}^{\prime} & b_{1,2}^{\prime} & \cdots & b_{1,\ell^{\prime}}^{\prime}\\
a_{2,1}^{\prime} & a_{2,2}^{\prime} & \cdots & a_{2,\ell}^{\prime} &
b_{2,1}^{\prime} & b_{2,2}^{\prime} & \cdots & b_{2,\ell^{\prime}}^{\prime}\\
\vdots & \vdots & \ddots & \vdots & \vdots & \vdots & \ddots & \vdots\\
a_{m,1}^{\prime} & a_{m,2}^{\prime} & \cdots & a_{m,\ell}^{\prime} &
b_{m,1}^{\prime} & b_{m,2}^{\prime} & \cdots & b_{m,\ell^{\prime}}^{\prime}\\
c_{1,1}^{\prime} & c_{1,2}^{\prime} & \cdots & c_{1,\ell}^{\prime} &
d_{1,1}^{\prime} & d_{1,2}^{\prime} & \cdots & d_{1,\ell^{\prime}}^{\prime}\\
c_{2,1}^{\prime} & c_{2,2}^{\prime} & \cdots & c_{2,\ell}^{\prime} &
d_{2,1}^{\prime} & d_{2,2}^{\prime} & \cdots & d_{2,\ell^{\prime}}^{\prime}\\
\vdots & \vdots & \ddots & \vdots & \vdots & \vdots & \ddots & \vdots\\
c_{m^{\prime},1}^{\prime} & c_{m^{\prime},2}^{\prime} & \cdots & c_{m^{\prime
},\ell}^{\prime} & d_{m^{\prime},1}^{\prime} & d_{m^{\prime},2}^{\prime} &
\cdots & d_{m^{\prime},\ell^{\prime}}^{\prime}%
\end{array}
\right) \\
&  =\left(
\begin{cases}
a_{i,j}^{\prime}, & \text{if }i\leq m\text{ \& }j\leq\ell;\\
b_{i,j-\ell}^{\prime}, & \text{if }i\leq m\text{ \& }j>\ell;\\
c_{i-m,j}^{\prime}, & \text{if }i>m\text{ \& }j\leq\ell;\\
d_{i-m,j-\ell}^{\prime}, & \text{if }i>m\text{ \& }j>\ell
\end{cases}
\right)  _{1\leq i\leq m+m^{\prime},\ 1\leq j\leq\ell+\ell^{\prime}}.
\end{align*}
Using these two equalities, we can compute the product $\left(
\begin{array}
[c]{cc}%
A & B\\
C & D
\end{array}
\right)  \left(
\begin{array}
[c]{cc}%
A^{\prime} & B^{\prime}\\
C^{\prime} & D^{\prime}%
\end{array}
\right)  $: Namely,%
\begin{align}
&  \left(
\begin{array}
[c]{cc}%
A & B\\
C & D
\end{array}
\right)  \left(
\begin{array}
[c]{cc}%
A^{\prime} & B^{\prime}\\
C^{\prime} & D^{\prime}%
\end{array}
\right) \nonumber\\
&  =\left(  \sum_{k=1}^{m+m^{\prime}}%
\begin{cases}
a_{i,k}, & \text{if }i\leq n\text{ \& }k\leq m;\\
b_{i,k-m}, & \text{if }i\leq n\text{ \& }k>m;\\
c_{i-n,k}, & \text{if }i>n\text{ \& }k\leq m;\\
d_{i-n,k-m}, & \text{if }i>n\text{ \& }k>m
\end{cases}%
\begin{cases}
a_{k,j}^{\prime}, & \text{if }k\leq m\text{ \& }j\leq\ell;\\
b_{k,j-\ell}^{\prime}, & \text{if }k\leq m\text{ \& }j>\ell;\\
c_{k-m,j}^{\prime}, & \text{if }k>m\text{ \& }j\leq\ell;\\
d_{k-m,j-\ell}^{\prime}, & \text{if }k>m\text{ \& }j>\ell
\end{cases}
\right)  _{1\leq i\leq n+n^{\prime},\ 1\leq j\leq\ell+\ell^{\prime}}.
\label{sol.block2x2.mult.short.3}%
\end{align}

On the other hand, we have $A=\left(  a_{i,j}\right)  _{1\leq i\leq n,\ 1\leq
j\leq m}$ and $A^{\prime}=\left(  a_{i,j}^{\prime}\right)  _{1\leq i\leq
m,\ 1\leq j\leq\ell}$. Hence, the definition of the product of two matrices
shows that%
\begin{equation}
AA^{\prime}=\left(  \sum_{k=1}^{m}a_{i,k}a_{k,j}^{\prime}\right)  _{1\leq
i\leq n,\ 1\leq j\leq\ell}. \label{sol.block2x2.mult.short.5a}%
\end{equation}

Also, we have $B=\left(  b_{i,j}\right)  _{1\leq i\leq n,\ 1\leq j\leq
m^{\prime}}$ and $C^{\prime}=\left(  c_{i,j}^{\prime}\right)  _{1\leq i\leq
m^{\prime},\ 1\leq j\leq\ell}$. Hence, the definition of the product of two
matrices shows that%
\begin{equation}
BC^{\prime}=\left(  \sum_{k=1}^{m^{\prime}}b_{i,k}c_{k,j}^{\prime}\right)
_{1\leq i\leq n,\ 1\leq j\leq\ell}. \label{sol.block2x2.mult.short.5b}%
\end{equation}

Adding the equalities (\ref{sol.block2x2.mult.short.5a}) and
(\ref{sol.block2x2.mult.short.5b}), we obtain%
\begin{align}
AA^{\prime}+BC^{\prime}  &  =\left(  \sum_{k=1}^{m}a_{i,k}a_{k,j}^{\prime
}\right)  _{1\leq i\leq n,\ 1\leq j\leq\ell}+\left(  \sum_{k=1}^{m^{\prime}%
}b_{i,k}c_{k,j}^{\prime}\right)  _{1\leq i\leq n,\ 1\leq j\leq\ell}\nonumber\\
&  =\left(  \sum_{k=1}^{m}a_{i,k}a_{k,j}^{\prime}+\underbrace{\sum
_{k=1}^{m^{\prime}}b_{i,k}c_{k,j}^{\prime}}_{\substack{=\sum_{k=m+1}%
^{m+m^{\prime}}b_{i,k-m}c_{k-m,j}^{\prime}\\\text{(here, we have substituted
}k-m\text{ for }k\text{ in the sum)}}}\right)  _{1\leq i\leq n,\ 1\leq
j\leq\ell}\nonumber\\
&  =\left(  \sum_{k=1}^{m}a_{i,k}a_{k,j}^{\prime}+\sum_{k=m+1}^{m+m^{\prime}%
}b_{i,k-m}c_{k-m,j}^{\prime}\right)  _{1\leq i\leq n,\ 1\leq j\leq\ell}.
\label{sol.block2x2.mult.short.5c}%
\end{align}

Similarly,%
\begin{equation}
AB^{\prime}+BD^{\prime}=\left(  \sum_{k=1}^{m}a_{i,k}b_{k,j}^{\prime}%
+\sum_{k=m+1}^{m+m^{\prime}}b_{i,k-m}d_{k-m,j}^{\prime}\right)  _{1\leq i\leq
n,\ 1\leq j\leq\ell^{\prime}}; \label{sol.block2x2.mult.short.6c}%
\end{equation}%
\begin{equation}
CA^{\prime}+DC^{\prime}=\left(  \sum_{k=1}^{m}c_{i,k}a_{k,j}^{\prime}%
+\sum_{k=m+1}^{m+m^{\prime}}d_{i,k-m}c_{k-m,j}^{\prime}\right)  _{1\leq i\leq
n^{\prime},\ 1\leq j\leq\ell}; \label{sol.block2x2.mult.short.7c}%
\end{equation}%
\begin{equation}
CB^{\prime}+DD^{\prime}=\left(  \sum_{k=1}^{m}c_{i,k}b_{k,j}^{\prime}%
+\sum_{k=m+1}^{m+m^{\prime}}d_{i,k-m}d_{k-m,j}^{\prime}\right)  _{1\leq i\leq
n^{\prime},\ 1\leq j\leq\ell^{\prime}}. \label{sol.block2x2.mult.short.8c}%
\end{equation}

Now, we have the four equalities (\ref{sol.block2x2.mult.short.5c}),
(\ref{sol.block2x2.mult.short.6c}), (\ref{sol.block2x2.mult.short.7c}) and
(\ref{sol.block2x2.mult.short.8c}). Hence, the definition of the block matrix
$\left(
\begin{array}
[c]{cc}%
AA^{\prime}+BC^{\prime} & AB^{\prime}+BD^{\prime}\\
CA^{\prime}+DC^{\prime} & CB^{\prime}+DD^{\prime}%
\end{array}
\right)  $ (or, more precisely, the equality (\ref{eq.def.block2x2.formal}),
applied to $n$, $n^{\prime}$, $\ell$, $\ell^{\prime}$, $AA^{\prime}%
+BC^{\prime}$, $AB^{\prime}+BD^{\prime}$, $CA^{\prime}+DC^{\prime}$,
$CB^{\prime}+DD^{\prime}$, $\sum_{k=1}^{m}a_{i,k}a_{k,j}^{\prime}+\sum
_{k=m+1}^{m+m^{\prime}}b_{i,k-m}c_{k-m,j}^{\prime}$, $\sum_{k=1}^{m}%
a_{i,k}b_{k,j}^{\prime}+\sum_{k=m+1}^{m+m^{\prime}}b_{i,k-m}d_{k-m,j}^{\prime
}$, \newline$\sum_{k=1}^{m}c_{i,k}a_{k,j}^{\prime}+\sum_{k=m+1}^{m+m^{\prime}%
}d_{i,k-m}c_{k-m,j}^{\prime}$ and $\sum_{k=1}^{m}c_{i,k}b_{k,j}^{\prime}%
+\sum_{k=m+1}^{m+m^{\prime}}d_{i,k-m}d_{k-m,j}^{\prime}$ instead of $n$,
$n^{\prime}$, $m$, $m^{\prime}$, $A$, $B$, $C$, $D$, $a_{i,j}$, $b_{i,j}$,
$c_{i,j}$ and $d_{i,j}$) shows that
\begin{align}
&  \left(
\begin{array}
[c]{cc}%
AA^{\prime}+BC^{\prime} & AB^{\prime}+BD^{\prime}\\
CA^{\prime}+DC^{\prime} & CB^{\prime}+DD^{\prime}%
\end{array}
\right) \nonumber\\
&  =\left(
\begin{cases}
\sum_{k=1}^{m}a_{i,k}a_{k,j}^{\prime}+\sum_{k=m+1}^{m+m^{\prime}}%
b_{i,k-m}c_{k-m,j}^{\prime}, & \text{if }i\leq n\text{ \& }j\leq\ell;\\
\sum_{k=1}^{m}a_{i,k}b_{k,j-\ell}^{\prime}+\sum_{k=m+1}^{m+m^{\prime}%
}b_{i,k-m}d_{k-m,j-\ell}^{\prime}, & \text{if }i\leq n\text{ \& }j>\ell;\\
\sum_{k=1}^{m}c_{i-n,k}a_{k,j}^{\prime}+\sum_{k=m+1}^{m+m^{\prime}}%
d_{i-n,k-m}c_{k-m,j}^{\prime}, & \text{if }i>n\text{ \& }j\leq\ell;\\
\sum_{k=1}^{m}c_{i-n,k}b_{k,j-\ell}^{\prime}+\sum_{k=m+1}^{m+m^{\prime}%
}d_{i-n,k-m}d_{k-m,j-\ell}^{\prime}, & \text{if }i>n\text{ \& }j>\ell
\end{cases}
\right)  _{1\leq i\leq n+n^{\prime},\ 1\leq j\leq\ell+\ell^{\prime}}.
\label{sol.block2x2.mult.short.9}%
\end{align}

But our goal is to prove that $\left(
\begin{array}
[c]{cc}%
A & B\\
C & D
\end{array}
\right)  \left(
\begin{array}
[c]{cc}%
A^{\prime} & B^{\prime}\\
C^{\prime} & D^{\prime}%
\end{array}
\right)  =\left(
\begin{array}
[c]{cc}%
AA^{\prime}+BC^{\prime} & AB^{\prime}+BD^{\prime}\\
CA^{\prime}+DC^{\prime} & CB^{\prime}+DD^{\prime}%
\end{array}
\right)  $. In other words, we want to prove that the left hand sides of the
equalities (\ref{sol.block2x2.mult.short.3}) and
(\ref{sol.block2x2.mult.short.9}) are equal. For this, it clearly suffices to
show that the right hand sides of these equalities are equal. In other words,
it suffices to show that every $\left(  i,j\right)  \in\left\{  1,2,\ldots
,n+n^{\prime}\right\}  \times\left\{  1,2,\ldots,\ell+\ell^{\prime}\right\}  $
satisfies%
\begin{align}
&  \sum_{k=1}^{m+m^{\prime}}%
\begin{cases}
a_{i,k}, & \text{if }i\leq n\text{ \& }k\leq m;\\
b_{i,k-m}, & \text{if }i\leq n\text{ \& }k>m;\\
c_{i-n,k}, & \text{if }i>n\text{ \& }k\leq m;\\
d_{i-n,k-m}, & \text{if }i>n\text{ \& }k>m
\end{cases}%
\begin{cases}
a_{k,j}^{\prime}, & \text{if }k\leq m\text{ \& }j\leq\ell;\\
b_{k,j-\ell}^{\prime}, & \text{if }k\leq m\text{ \& }j>\ell;\\
c_{k-m,j}^{\prime}, & \text{if }k>m\text{ \& }j\leq\ell;\\
d_{k-m,j-\ell}^{\prime}, & \text{if }k>m\text{ \& }j>\ell
\end{cases}
\nonumber\\
&  =%
\begin{cases}
\sum_{k=1}^{m}a_{i,k}a_{k,j}^{\prime}+\sum_{k=m+1}^{m+m^{\prime}}%
b_{i,k-m}c_{k-m,j}^{\prime}, & \text{if }i\leq n\text{ \& }j\leq\ell;\\
\sum_{k=1}^{m}a_{i,k}b_{k,j-\ell}^{\prime}+\sum_{k=m+1}^{m+m^{\prime}%
}b_{i,k-m}d_{k-m,j-\ell}^{\prime}, & \text{if }i\leq n\text{ \& }j>\ell;\\
\sum_{k=1}^{m}c_{i-n,k}a_{k,j}^{\prime}+\sum_{k=m+1}^{m+m^{\prime}}%
d_{i-n,k-m}c_{k-m,j}^{\prime}, & \text{if }i>n\text{ \& }j\leq\ell;\\
\sum_{k=1}^{m}c_{i-n,k}b_{k,j-\ell}^{\prime}+\sum_{k=m+1}^{m+m^{\prime}%
}d_{i-n,k-m}d_{k-m,j-\ell}^{\prime}, & \text{if }i>n\text{ \& }j>\ell
\end{cases}
. \label{sol.block2x2.mult.short.entrywise}%
\end{align}

[\textit{Proof of (\ref{sol.block2x2.mult.short.entrywise}):} Let $\left(
i,j\right)  \in\left\{  1,2,\ldots,n+n^{\prime}\right\}  \times\left\{
1,2,\ldots,\ell+\ell^{\prime}\right\}  $. Thus, \newline$i\in\left\{
1,2,\ldots,n+n^{\prime}\right\}  $ and $j\in\left\{  1,2,\ldots,\ell
+\ell^{\prime}\right\}  $. We must be in one of the following four cases:

\textit{Case 1:} We have $i\leq n$ and $j\leq\ell$.

\textit{Case 2:} We have $i\leq n$ and $j>\ell$.

\textit{Case 3:} We have $i>n$ and $j\leq\ell$.

\textit{Case 4:} We have $i>n$ and $j>\ell$.

All four cases are completely analogous; we thus will only show how to deal
with Case 1. In this case, we have $i\leq n$ and $j\leq\ell$. Now, comparing%
\begin{align*}
&  \sum_{k=1}^{m+m^{\prime}}%
\begin{cases}
a_{i,k}, & \text{if }i\leq n\text{ \& }k\leq m;\\
b_{i,k-m}, & \text{if }i\leq n\text{ \& }k>m;\\
c_{i-n,k}, & \text{if }i>n\text{ \& }k\leq m;\\
d_{i-n,k-m}, & \text{if }i>n\text{ \& }k>m
\end{cases}%
\begin{cases}
a_{k,j}^{\prime}, & \text{if }k\leq m\text{ \& }j\leq\ell;\\
b_{k,j-\ell}^{\prime}, & \text{if }k\leq m\text{ \& }j>\ell;\\
c_{k-m,j}^{\prime}, & \text{if }k>m\text{ \& }j\leq\ell;\\
d_{k-m,j-\ell}^{\prime}, & \text{if }k>m\text{ \& }j>\ell
\end{cases}
\\
&  =\sum_{k=1}^{m}\underbrace{%
\begin{cases}
a_{i,k}, & \text{if }i\leq n\text{ \& }k\leq m;\\
b_{i,k-m}, & \text{if }i\leq n\text{ \& }k>m;\\
c_{i-n,k}, & \text{if }i>n\text{ \& }k\leq m;\\
d_{i-n,k-m}, & \text{if }i>n\text{ \& }k>m
\end{cases}
}_{\substack{=a_{i,k}\\\text{(since }i\leq n\text{ and }k\leq m\text{)}%
}}\underbrace{%
\begin{cases}
a_{k,j}^{\prime}, & \text{if }k\leq m\text{ \& }j\leq\ell;\\
b_{k,j-\ell}^{\prime}, & \text{if }k\leq m\text{ \& }j>\ell;\\
c_{k-m,j}^{\prime}, & \text{if }k>m\text{ \& }j\leq\ell;\\
d_{k-m,j-\ell}^{\prime}, & \text{if }k>m\text{ \& }j>\ell
\end{cases}
}_{\substack{=a_{k,j}^{\prime}\\\text{(since }k\leq m\text{ and }j\leq
\ell\text{)}}}\\
&  \ \ \ \ \ \ \ \ \ \ +\sum_{k=m+1}^{m+m^{\prime}}\underbrace{%
\begin{cases}
a_{i,k}, & \text{if }i\leq n\text{ \& }k\leq m;\\
b_{i,k-m}, & \text{if }i\leq n\text{ \& }k>m;\\
c_{i-n,k}, & \text{if }i>n\text{ \& }k\leq m;\\
d_{i-n,k-m}, & \text{if }i>n\text{ \& }k>m
\end{cases}
}_{\substack{=b_{i,k-m}\\\text{(since }i\leq n\text{ and }k>m\text{)}%
}}\underbrace{%
\begin{cases}
a_{k,j}^{\prime}, & \text{if }k\leq m\text{ \& }j\leq\ell;\\
b_{k,j-\ell}^{\prime}, & \text{if }k\leq m\text{ \& }j>\ell;\\
c_{k-m,j}^{\prime}, & \text{if }k>m\text{ \& }j\leq\ell;\\
d_{k-m,j-\ell}^{\prime}, & \text{if }k>m\text{ \& }j>\ell
\end{cases}
}_{\substack{=c_{k-m,j}^{\prime}\\\text{(since }k>m\text{ and }j\leq
\ell\text{)}}}\\
&  \ \ \ \ \ \ \ \ \ \ \left(  \text{since }0\leq m\leq m+m^{\prime}\right) \\
&  =\sum_{k=1}^{m}a_{i,k}a_{k,j}^{\prime}+\sum_{k=m+1}^{m+m^{\prime}}%
b_{i,k-m}c_{k-m,j}^{\prime}%
\end{align*}
with%
\begin{align*}
&
\begin{cases}
\sum_{k=1}^{m}a_{i,k}a_{k,j}^{\prime}+\sum_{k=m+1}^{m+m^{\prime}}%
b_{i,k-m}c_{k-m,j}^{\prime}, & \text{if }i\leq n\text{ \& }j\leq\ell;\\
\sum_{k=1}^{m}a_{i,k}b_{k,j-\ell}^{\prime}+\sum_{k=m+1}^{m+m^{\prime}%
}b_{i,k-m}d_{k-m,j-\ell}^{\prime}, & \text{if }i\leq n\text{ \& }j>\ell;\\
\sum_{k=1}^{m}c_{i-n,k}a_{k,j}^{\prime}+\sum_{k=m+1}^{m+m^{\prime}}%
d_{i-n,k-m}c_{k-m,j}^{\prime}, & \text{if }i>n\text{ \& }j\leq\ell;\\
\sum_{k=1}^{m}c_{i-n,k}b_{k,j-\ell}^{\prime}+\sum_{k=m+1}^{m+m^{\prime}%
}d_{i-n,k-m}d_{k-m,j-\ell}^{\prime}, & \text{if }i>n\text{ \& }j>\ell
\end{cases}
\\
&  =\sum_{k=1}^{m}a_{i,k}a_{k,j}^{\prime}+\sum_{k=m+1}^{m+m^{\prime}}%
b_{i,k-m}c_{k-m,j}^{\prime}\ \ \ \ \ \ \ \ \ \ \left(  \text{since }i\leq
n\text{ and }j\leq\ell\right)  ,
\end{align*}
we obtain precisely (\ref{sol.block2x2.mult.short.entrywise}). Thus,
(\ref{sol.block2x2.mult.short.entrywise}) is proven in Case 1. As I said, the
other three cases are similar, and so (\ref{sol.block2x2.mult.short.entrywise}%
) is proven.]

From (\ref{sol.block2x2.mult.short.entrywise}), we see that the right hand
sides of the equalities (\ref{sol.block2x2.mult.short.3}) and
(\ref{sol.block2x2.mult.short.9}) are equal. Hence, so are their left hand
sides. In other words, $\left(
\begin{array}
[c]{cc}%
A & B\\
C & D
\end{array}
\right)  \left(
\begin{array}
[c]{cc}%
A^{\prime} & B^{\prime}\\
C^{\prime} & D^{\prime}%
\end{array}
\right)  =\left(
\begin{array}
[c]{cc}%
AA^{\prime}+BC^{\prime} & AB^{\prime}+BD^{\prime}\\
CA^{\prime}+DC^{\prime} & CB^{\prime}+DD^{\prime}%
\end{array}
\right)  $. This solves Exercise \ref{exe.block2x2.mult}.
\end{proof}
\end{vershort}

\begin{verlong}
\begin{proof}
[Solution to Exercise \ref{exe.block2x2.mult}.]Write the $n\times m$-matrix
$A$ in the form $A=\left(  a_{i,j}\right)  _{1\leq i\leq n,\ 1\leq j\leq m}$.

Write the $n\times m^{\prime}$-matrix $B$ in the form $B=\left(
b_{i,j}\right)  _{1\leq i\leq n,\ 1\leq j\leq m^{\prime}}$.

Write the $n^{\prime}\times m$-matrix $C$ in the form $C=\left(
c_{i,j}\right)  _{1\leq i\leq n^{\prime},\ 1\leq j\leq m}$.

Write the $n^{\prime}\times m^{\prime}$-matrix $D$ in the form $D=\left(
d_{i,j}\right)  _{1\leq i\leq n^{\prime},\ 1\leq j\leq m^{\prime}}$.

Write the $m\times\ell$-matrix $A^{\prime}$ in the form $A^{\prime}=\left(
a_{i,j}^{\prime}\right)  _{1\leq i\leq m,\ 1\leq j\leq\ell}$.

Write the $m\times\ell^{\prime}$-matrix $B^{\prime}$ in the form $B^{\prime
}=\left(  b_{i,j}^{\prime}\right)  _{1\leq i\leq m,\ 1\leq j\leq\ell^{\prime}%
}$.

Write the $m^{\prime}\times\ell$-matrix $C^{\prime}$ in the form $C^{\prime
}=\left(  c_{i,j}^{\prime}\right)  _{1\leq i\leq m^{\prime},\ 1\leq j\leq\ell
}$.

Write the $m^{\prime}\times\ell^{\prime}$-matrix $D^{\prime}$ in the form
$D^{\prime}=\left(  d_{i,j}^{\prime}\right)  _{1\leq i\leq m^{\prime},\ 1\leq
j\leq\ell^{\prime}}$.

We have $A^{\prime}=\left(  a_{i,j}^{\prime}\right)  _{1\leq i\leq m,\ 1\leq
j\leq\ell}$, $B^{\prime}=\left(  b_{i,j}^{\prime}\right)  _{1\leq i\leq
m,\ 1\leq j\leq\ell^{\prime}}$, $C^{\prime}=\left(  c_{i,j}^{\prime}\right)
_{1\leq i\leq m^{\prime},\ 1\leq j\leq\ell}$ and $D^{\prime}=\left(
d_{i,j}^{\prime}\right)  _{1\leq i\leq m^{\prime},\ 1\leq j\leq\ell^{\prime}}%
$. Thus, (\ref{eq.def.block2x2.formal}) (applied to $m$, $m^{\prime}$, $\ell$,
$\ell^{\prime}$, $A^{\prime}$, $B^{\prime}$, $C^{\prime}$, $D^{\prime}$,
$a_{i,j}^{\prime}$, $b_{i,j}^{\prime}$, $c_{i,j}^{\prime}$ and $d_{i,j}%
^{\prime}$ instead of $n$, $n^{\prime}$, $m$, $m^{\prime}$, $A$, $B$, $C$,
$D$, $a_{i,j}$, $b_{i,j}$, $c_{i,j}$ and $d_{i,j}$) shows that
\begin{equation}
\left(
\begin{array}
[c]{cc}%
A^{\prime} & B^{\prime}\\
C^{\prime} & D^{\prime}%
\end{array}
\right)  =\left(
\begin{cases}
a_{i,j}^{\prime}, & \text{if }i\leq m\text{ and }j\leq\ell;\\
b_{i,j-\ell}^{\prime}, & \text{if }i\leq m\text{ and }j>\ell;\\
c_{i-m,j}^{\prime}, & \text{if }i>m\text{ and }j\leq\ell;\\
d_{i-m,j-\ell}^{\prime}, & \text{if }i>m\text{ and }j>\ell
\end{cases}
\right)  _{1\leq i\leq m+m^{\prime},\ 1\leq j\leq\ell+\ell^{\prime}}.
\label{sol.block2x2.mult.1}%
\end{equation}

Now, we have (\ref{eq.def.block2x2.formal}) and (\ref{sol.block2x2.mult.1}).
Thus, the definition of the product of two matrices shows that%
\begin{align}
&  \left(
\begin{array}
[c]{cc}%
A & B\\
C & D
\end{array}
\right)  \left(
\begin{array}
[c]{cc}%
A^{\prime} & B^{\prime}\\
C^{\prime} & D^{\prime}%
\end{array}
\right) \nonumber\\
&  =\left(  \sum_{k=1}^{m+m^{\prime}}%
\begin{cases}
a_{i,k}, & \text{if }i\leq n\text{ and }k\leq m;\\
b_{i,k-m}, & \text{if }i\leq n\text{ and }k>m;\\
c_{i-n,k}, & \text{if }i>n\text{ and }k\leq m;\\
d_{i-n,k-m}, & \text{if }i>n\text{ and }k>m
\end{cases}%
\begin{cases}
a_{k,j}^{\prime}, & \text{if }k\leq m\text{ and }j\leq\ell;\\
b_{k,j-\ell}^{\prime}, & \text{if }k\leq m\text{ and }j>\ell;\\
c_{k-m,j}^{\prime}, & \text{if }k>m\text{ and }j\leq\ell;\\
d_{k-m,j-\ell}^{\prime}, & \text{if }k>m\text{ and }j>\ell
\end{cases}
\right)  _{1\leq i\leq n+n^{\prime},\ 1\leq j\leq\ell+\ell^{\prime}}.
\label{sol.block2x2.mult.3}%
\end{align}

On the other hand, we have $A=\left(  a_{i,j}\right)  _{1\leq i\leq n,\ 1\leq
j\leq m}$ and $A^{\prime}=\left(  a_{i,j}^{\prime}\right)  _{1\leq i\leq
m,\ 1\leq j\leq\ell}$. Hence, the definition of the product of two matrices
shows that%
\begin{equation}
AA^{\prime}=\left(  \sum_{k=1}^{m}a_{i,k}a_{k,j}^{\prime}\right)  _{1\leq
i\leq n,\ 1\leq j\leq\ell}. \label{sol.block2x2.mult.5a}%
\end{equation}

Also, we have $B=\left(  b_{i,j}\right)  _{1\leq i\leq n,\ 1\leq j\leq
m^{\prime}}$ and $C^{\prime}=\left(  c_{i,j}^{\prime}\right)  _{1\leq i\leq
m^{\prime},\ 1\leq j\leq\ell}$. Hence, the definition of the product of two
matrices shows that%
\begin{equation}
BC^{\prime}=\left(  \sum_{k=1}^{m^{\prime}}b_{i,k}c_{k,j}^{\prime}\right)
_{1\leq i\leq n,\ 1\leq j\leq\ell}. \label{sol.block2x2.mult.5b}%
\end{equation}

Adding the equalities (\ref{sol.block2x2.mult.5a}) and
(\ref{sol.block2x2.mult.5b}), we obtain%
\begin{align}
AA^{\prime}+BC^{\prime}  &  =\left(  \sum_{k=1}^{m}a_{i,k}a_{k,j}^{\prime
}\right)  _{1\leq i\leq n,\ 1\leq j\leq\ell}+\left(  \sum_{k=1}^{m^{\prime}%
}b_{i,k}c_{k,j}^{\prime}\right)  _{1\leq i\leq n,\ 1\leq j\leq\ell}\nonumber\\
&  =\left(  \sum_{k=1}^{m}a_{i,k}a_{k,j}^{\prime}+\underbrace{\sum
_{k=1}^{m^{\prime}}b_{i,k}c_{k,j}^{\prime}}_{\substack{=\sum_{k=m+1}%
^{m+m^{\prime}}b_{i,k-m}c_{k-m,j}^{\prime}\\\text{(here, we have substituted
}k-m\text{ for }k\text{ in the sum)}}}\right)  _{1\leq i\leq n,\ 1\leq
j\leq\ell}\nonumber\\
&  =\left(  \sum_{k=1}^{m}a_{i,k}a_{k,j}^{\prime}+\sum_{k=m+1}^{m+m^{\prime}%
}b_{i,k-m}c_{k-m,j}^{\prime}\right)  _{1\leq i\leq n,\ 1\leq j\leq\ell}.
\label{sol.block2x2.mult.5c}%
\end{align}

Furthermore, we have $A=\left(  a_{i,j}\right)  _{1\leq i\leq n,\ 1\leq j\leq
m}$ and $B^{\prime}=\left(  b_{i,j}^{\prime}\right)  _{1\leq i\leq m,\ 1\leq
j\leq\ell^{\prime}}$. Hence, the definition of the product of two matrices
shows that%
\begin{equation}
AB^{\prime}=\left(  \sum_{k=1}^{m}a_{i,k}b_{k,j}^{\prime}\right)  _{1\leq
i\leq n,\ 1\leq j\leq\ell^{\prime}}. \label{sol.block2x2.mult.6a}%
\end{equation}

Also, we have $B=\left(  b_{i,j}\right)  _{1\leq i\leq n,\ 1\leq j\leq
m^{\prime}}$ and $D^{\prime}=\left(  d_{i,j}^{\prime}\right)  _{1\leq i\leq
m^{\prime},\ 1\leq j\leq\ell^{\prime}}$. Hence, the definition of the product
of two matrices shows that%
\begin{equation}
BD^{\prime}=\left(  \sum_{k=1}^{m^{\prime}}b_{i,k}d_{k,j}^{\prime}\right)
_{1\leq i\leq n,\ 1\leq j\leq\ell^{\prime}}. \label{sol.block2x2.mult.6b}%
\end{equation}

Adding the equalities (\ref{sol.block2x2.mult.6a}) and
(\ref{sol.block2x2.mult.6b}), we obtain%
\begin{align}
AB^{\prime}+BD^{\prime}  &  =\left(  \sum_{k=1}^{m}a_{i,k}b_{k,j}^{\prime
}\right)  _{1\leq i\leq n,\ 1\leq j\leq\ell^{\prime}}+\left(  \sum
_{k=1}^{m^{\prime}}b_{i,k}d_{k,j}^{\prime}\right)  _{1\leq i\leq n,\ 1\leq
j\leq\ell^{\prime}}\nonumber\\
&  =\left(  \sum_{k=1}^{m}a_{i,k}b_{k,j}^{\prime}+\underbrace{\sum
_{k=1}^{m^{\prime}}b_{i,k}d_{k,j}^{\prime}}_{\substack{=\sum_{k=m+1}%
^{m+m^{\prime}}b_{i,k-m}d_{k-m,j}^{\prime}\\\text{(here, we have substituted
}k-m\text{ for }k\text{ in the sum)}}}\right)  _{1\leq i\leq n,\ 1\leq
j\leq\ell^{\prime}}\nonumber\\
&  =\left(  \sum_{k=1}^{m}a_{i,k}b_{k,j}^{\prime}+\sum_{k=m+1}^{m+m^{\prime}%
}b_{i,k-m}d_{k-m,j}^{\prime}\right)  _{1\leq i\leq n,\ 1\leq j\leq\ell
^{\prime}}. \label{sol.block2x2.mult.6c}%
\end{align}

Furthermore, we have $C=\left(  c_{i,j}\right)  _{1\leq i\leq n^{\prime
},\ 1\leq j\leq m}$ and $A^{\prime}=\left(  a_{i,j}^{\prime}\right)  _{1\leq
i\leq m,\ 1\leq j\leq\ell}$. Hence, the definition of the product of two
matrices shows that%
\begin{equation}
CA^{\prime}=\left(  \sum_{k=1}^{m}c_{i,k}a_{k,j}^{\prime}\right)  _{1\leq
i\leq n^{\prime},\ 1\leq j\leq\ell}. \label{sol.block2x2.mult.7a}%
\end{equation}

Also, we have $D=\left(  d_{i,j}\right)  _{1\leq i\leq n^{\prime},\ 1\leq
j\leq m^{\prime}}$ and $C^{\prime}=\left(  c_{i,j}^{\prime}\right)  _{1\leq
i\leq m^{\prime},\ 1\leq j\leq\ell}$. Hence, the definition of the product of
two matrices shows that%
\begin{equation}
DC^{\prime}=\left(  \sum_{k=1}^{m^{\prime}}d_{i,k}c_{k,j}^{\prime}\right)
_{1\leq i\leq n^{\prime},\ 1\leq j\leq\ell}. \label{sol.block2x2.mult.7b}%
\end{equation}

Adding the equalities (\ref{sol.block2x2.mult.7a}) and
(\ref{sol.block2x2.mult.7b}), we obtain%
\begin{align}
CA^{\prime}+DC^{\prime}  &  =\left(  \sum_{k=1}^{m}c_{i,k}a_{k,j}^{\prime
}\right)  _{1\leq i\leq n^{\prime},\ 1\leq j\leq\ell}+\left(  \sum
_{k=1}^{m^{\prime}}d_{i,k}c_{k,j}^{\prime}\right)  _{1\leq i\leq n^{\prime
},\ 1\leq j\leq\ell}\nonumber\\
&  =\left(  \sum_{k=1}^{m}c_{i,k}a_{k,j}^{\prime}+\underbrace{\sum
_{k=1}^{m^{\prime}}d_{i,k}c_{k,j}^{\prime}}_{\substack{=\sum_{k=m+1}%
^{m+m^{\prime}}d_{i,k-m}c_{k-m,j}^{\prime}\\\text{(here, we have substituted
}k-m\text{ for }k\text{ in the sum)}}}\right)  _{1\leq i\leq n^{\prime
},\ 1\leq j\leq\ell}\nonumber\\
&  =\left(  \sum_{k=1}^{m}c_{i,k}a_{k,j}^{\prime}+\sum_{k=m+1}^{m+m^{\prime}%
}d_{i,k-m}c_{k-m,j}^{\prime}\right)  _{1\leq i\leq n^{\prime},\ 1\leq
j\leq\ell}. \label{sol.block2x2.mult.7c}%
\end{align}

Furthermore, we have $C=\left(  c_{i,j}\right)  _{1\leq i\leq n^{\prime
},\ 1\leq j\leq m}$ and $B^{\prime}=\left(  b_{i,j}^{\prime}\right)  _{1\leq
i\leq m,\ 1\leq j\leq\ell^{\prime}}$. Hence, the definition of the product of
two matrices shows that%
\begin{equation}
CB^{\prime}=\left(  \sum_{k=1}^{m}c_{i,k}b_{k,j}^{\prime}\right)  _{1\leq
i\leq n^{\prime},\ 1\leq j\leq\ell^{\prime}}. \label{sol.block2x2.mult.8a}%
\end{equation}

Also, we have $D=\left(  d_{i,j}\right)  _{1\leq i\leq n^{\prime},\ 1\leq
j\leq m^{\prime}}$ and $D^{\prime}=\left(  d_{i,j}^{\prime}\right)  _{1\leq
i\leq m^{\prime},\ 1\leq j\leq\ell^{\prime}}$. Hence, the definition of the
product of two matrices shows that%
\begin{equation}
DD^{\prime}=\left(  \sum_{k=1}^{m^{\prime}}d_{i,k}d_{k,j}^{\prime}\right)
_{1\leq i\leq n^{\prime},\ 1\leq j\leq\ell^{\prime}}.
\label{sol.block2x2.mult.8b}%
\end{equation}

Adding the equalities (\ref{sol.block2x2.mult.8a}) and
(\ref{sol.block2x2.mult.8b}), we obtain%
\begin{align}
CB^{\prime}+DD^{\prime}  &  =\left(  \sum_{k=1}^{m}c_{i,k}b_{k,j}^{\prime
}\right)  _{1\leq i\leq n^{\prime},\ 1\leq j\leq\ell^{\prime}}+\left(
\sum_{k=1}^{m^{\prime}}d_{i,k}d_{k,j}^{\prime}\right)  _{1\leq i\leq
n^{\prime},\ 1\leq j\leq\ell^{\prime}}\nonumber\\
&  =\left(  \sum_{k=1}^{m}c_{i,k}b_{k,j}^{\prime}+\underbrace{\sum
_{k=1}^{m^{\prime}}d_{i,k}d_{k,j}^{\prime}}_{\substack{=\sum_{k=m+1}%
^{m+m^{\prime}}d_{i,k-m}d_{k-m,j}^{\prime}\\\text{(here, we have substituted
}k-m\text{ for }k\text{ in the sum)}}}\right)  _{1\leq i\leq n^{\prime
},\ 1\leq j\leq\ell^{\prime}}\nonumber\\
&  =\left(  \sum_{k=1}^{m}c_{i,k}b_{k,j}^{\prime}+\sum_{k=m+1}^{m+m^{\prime}%
}d_{i,k-m}d_{k-m,j}^{\prime}\right)  _{1\leq i\leq n^{\prime},\ 1\leq
j\leq\ell^{\prime}}. \label{sol.block2x2.mult.8c}%
\end{align}

Now, we have the four equalities (\ref{sol.block2x2.mult.5c}),
(\ref{sol.block2x2.mult.6c}), (\ref{sol.block2x2.mult.7c}) and
(\ref{sol.block2x2.mult.8c}). Hence, (\ref{eq.def.block2x2.formal}) (applied
to $n$, $n^{\prime}$, $\ell$, $\ell^{\prime}$, $AA^{\prime}+BC^{\prime}$,
$AB^{\prime}+BD^{\prime}$, $CA^{\prime}+DC^{\prime}$, $CB^{\prime}+DD^{\prime
}$, $\sum_{k=1}^{m}a_{i,k}a_{k,j}^{\prime}+\sum_{k=m+1}^{m+m^{\prime}%
}b_{i,k-m}c_{k-m,j}^{\prime}$, $\sum_{k=1}^{m}a_{i,k}b_{k,j}^{\prime}%
+\sum_{k=m+1}^{m+m^{\prime}}b_{i,k-m}d_{k-m,j}^{\prime}$, $\sum_{k=1}%
^{m}c_{i,k}a_{k,j}^{\prime}+\sum_{k=m+1}^{m+m^{\prime}}d_{i,k-m}%
c_{k-m,j}^{\prime}$ and $\sum_{k=1}^{m}c_{i,k}b_{k,j}^{\prime}+\sum
_{k=m+1}^{m+m^{\prime}}d_{i,k-m}d_{k-m,j}^{\prime}$ instead of $n$,
$n^{\prime}$, $m$, $m^{\prime}$, $A$, $B$, $C$, $D$, $a_{i,j}$, $b_{i,j}$,
$c_{i,j}$ and $d_{i,j}$) shows that
\begin{align}
&  \left(
\begin{array}
[c]{cc}%
AA^{\prime}+BC^{\prime} & AB^{\prime}+BD^{\prime}\\
CA^{\prime}+DC^{\prime} & CB^{\prime}+DD^{\prime}%
\end{array}
\right) \nonumber\\
&  =\left(
\begin{cases}
\sum_{k=1}^{m}a_{i,k}a_{k,j}^{\prime}+\sum_{k=m+1}^{m+m^{\prime}}%
b_{i,k-m}c_{k-m,j}^{\prime}, & \text{if }i\leq n\text{ and }j\leq\ell;\\
\sum_{k=1}^{m}a_{i,k}b_{k,j-\ell}^{\prime}+\sum_{k=m+1}^{m+m^{\prime}%
}b_{i,k-m}d_{k-m,j-\ell}^{\prime}, & \text{if }i\leq n\text{ and }j>\ell;\\
\sum_{k=1}^{m}c_{i-n,k}a_{k,j}^{\prime}+\sum_{k=m+1}^{m+m^{\prime}}%
d_{i-n,k-m}c_{k-m,j}^{\prime}, & \text{if }i>n\text{ and }j\leq\ell;\\
\sum_{k=1}^{m}c_{i-n,k}b_{k,j-\ell}^{\prime}+\sum_{k=m+1}^{m+m^{\prime}%
}d_{i-n,k-m}d_{k-m,j-\ell}^{\prime}, & \text{if }i>n\text{ and }j>\ell
\end{cases}
\right)  _{1\leq i\leq n+n^{\prime},\ 1\leq j\leq\ell+\ell^{\prime}}.
\label{sol.block2x2.mult.9}%
\end{align}

Now, we shall show that every $\left(  i,j\right)  \in\left\{  1,2,\ldots
,n+n^{\prime}\right\}  \times\left\{  1,2,\ldots,\ell+\ell^{\prime}\right\}  $
satisfies%
\begin{align}
&  \sum_{k=1}^{m+m^{\prime}}%
\begin{cases}
a_{i,k}, & \text{if }i\leq n\text{ and }k\leq m;\\
b_{i,k-m}, & \text{if }i\leq n\text{ and }k>m;\\
c_{i-n,k}, & \text{if }i>n\text{ and }k\leq m;\\
d_{i-n,k-m}, & \text{if }i>n\text{ and }k>m
\end{cases}%
\begin{cases}
a_{k,j}^{\prime}, & \text{if }k\leq m\text{ and }j\leq\ell;\\
b_{k,j-\ell}^{\prime}, & \text{if }k\leq m\text{ and }j>\ell;\\
c_{k-m,j}^{\prime}, & \text{if }k>m\text{ and }j\leq\ell;\\
d_{k-m,j-\ell}^{\prime}, & \text{if }k>m\text{ and }j>\ell
\end{cases}
\nonumber\\
&  =%
\begin{cases}
\sum_{k=1}^{m}a_{i,k}a_{k,j}^{\prime}+\sum_{k=m+1}^{m+m^{\prime}}%
b_{i,k-m}c_{k-m,j}^{\prime}, & \text{if }i\leq n\text{ and }j\leq\ell;\\
\sum_{k=1}^{m}a_{i,k}b_{k,j-\ell}^{\prime}+\sum_{k=m+1}^{m+m^{\prime}%
}b_{i,k-m}d_{k-m,j-\ell}^{\prime}, & \text{if }i\leq n\text{ and }j>\ell;\\
\sum_{k=1}^{m}c_{i-n,k}a_{k,j}^{\prime}+\sum_{k=m+1}^{m+m^{\prime}}%
d_{i-n,k-m}c_{k-m,j}^{\prime}, & \text{if }i>n\text{ and }j\leq\ell;\\
\sum_{k=1}^{m}c_{i-n,k}b_{k,j-\ell}^{\prime}+\sum_{k=m+1}^{m+m^{\prime}%
}d_{i-n,k-m}d_{k-m,j-\ell}^{\prime}, & \text{if }i>n\text{ and }j>\ell
\end{cases}
. \label{sol.block2x2.mult.entrywise}%
\end{align}

[\textit{Proof of (\ref{sol.block2x2.mult.entrywise}):} Let $\left(
i,j\right)  \in\left\{  1,2,\ldots,n+n^{\prime}\right\}  \times\left\{
1,2,\ldots,\ell+\ell^{\prime}\right\}  $. Thus, $i\in\left\{  1,2,\ldots
,n+n^{\prime}\right\}  $ and $j\in\left\{  1,2,\ldots,\ell+\ell^{\prime
}\right\}  $. We must be in one of the following two cases:

\textit{Case 1:} We have $i\leq n$.

\textit{Case 2:} We have $i>n$.

Let us first consider Case 1. In this case, we have $i\leq n$. Now, we must be
in one of the following two subcases:

\textit{Subcase 1.1:} We have $j\leq\ell$.

\textit{Subcase 1.2:} We have $j>\ell$.

Let us first consider Subcase 1.1. In this Subcase, we have $j\leq\ell$. Now,
comparing%
\begin{align*}
&  \sum_{k=1}^{m+m^{\prime}}%
\begin{cases}
a_{i,k}, & \text{if }i\leq n\text{ and }k\leq m;\\
b_{i,k-m}, & \text{if }i\leq n\text{ and }k>m;\\
c_{i-n,k}, & \text{if }i>n\text{ and }k\leq m;\\
d_{i-n,k-m}, & \text{if }i>n\text{ and }k>m
\end{cases}%
\begin{cases}
a_{k,j}^{\prime}, & \text{if }k\leq m\text{ and }j\leq\ell;\\
b_{k,j-\ell}^{\prime}, & \text{if }k\leq m\text{ and }j>\ell;\\
c_{k-m,j}^{\prime}, & \text{if }k>m\text{ and }j\leq\ell;\\
d_{k-m,j-\ell}^{\prime}, & \text{if }k>m\text{ and }j>\ell
\end{cases}
\\
&  =\sum_{k=1}^{m}\underbrace{%
\begin{cases}
a_{i,k}, & \text{if }i\leq n\text{ and }k\leq m;\\
b_{i,k-m}, & \text{if }i\leq n\text{ and }k>m;\\
c_{i-n,k}, & \text{if }i>n\text{ and }k\leq m;\\
d_{i-n,k-m}, & \text{if }i>n\text{ and }k>m
\end{cases}
}_{\substack{=a_{i,k}\\\text{(since }i\leq n\text{ and }k\leq m\text{)}%
}}\underbrace{%
\begin{cases}
a_{k,j}^{\prime}, & \text{if }k\leq m\text{ and }j\leq\ell;\\
b_{k,j-\ell}^{\prime}, & \text{if }k\leq m\text{ and }j>\ell;\\
c_{k-m,j}^{\prime}, & \text{if }k>m\text{ and }j\leq\ell;\\
d_{k-m,j-\ell}^{\prime}, & \text{if }k>m\text{ and }j>\ell
\end{cases}
}_{\substack{=a_{k,j}^{\prime}\\\text{(since }k\leq m\text{ and }j\leq
\ell\text{)}}}\\
&  \ \ \ \ \ \ \ \ \ \ +\sum_{k=m+1}^{m+m^{\prime}}\underbrace{%
\begin{cases}
a_{i,k}, & \text{if }i\leq n\text{ and }k\leq m;\\
b_{i,k-m}, & \text{if }i\leq n\text{ and }k>m;\\
c_{i-n,k}, & \text{if }i>n\text{ and }k\leq m;\\
d_{i-n,k-m}, & \text{if }i>n\text{ and }k>m
\end{cases}
}_{\substack{=b_{i,k-m}\\\text{(since }i\leq n\text{ and }k>m\text{)}%
}}\underbrace{%
\begin{cases}
a_{k,j}^{\prime}, & \text{if }k\leq m\text{ and }j\leq\ell;\\
b_{k,j-\ell}^{\prime}, & \text{if }k\leq m\text{ and }j>\ell;\\
c_{k-m,j}^{\prime}, & \text{if }k>m\text{ and }j\leq\ell;\\
d_{k-m,j-\ell}^{\prime}, & \text{if }k>m\text{ and }j>\ell
\end{cases}
}_{\substack{=c_{k-m,j}^{\prime}\\\text{(since }k>m\text{ and }j\leq
\ell\text{)}}}\\
&  \ \ \ \ \ \ \ \ \ \ \left(  \text{since }0\leq m\leq m+m^{\prime}\right) \\
&  =\sum_{k=1}^{m}a_{i,k}a_{k,j}^{\prime}+\sum_{k=m+1}^{m+m^{\prime}}%
b_{i,k-m}c_{k-m,j}^{\prime}%
\end{align*}
with%
\begin{align*}
&
\begin{cases}
\sum_{k=1}^{m}a_{i,k}a_{k,j}^{\prime}+\sum_{k=m+1}^{m+m^{\prime}}%
b_{i,k-m}c_{k-m,j}^{\prime}, & \text{if }i\leq n\text{ and }j\leq\ell;\\
\sum_{k=1}^{m}a_{i,k}b_{k,j-\ell}^{\prime}+\sum_{k=m+1}^{m+m^{\prime}%
}b_{i,k-m}d_{k-m,j-\ell}^{\prime}, & \text{if }i\leq n\text{ and }j>\ell;\\
\sum_{k=1}^{m}c_{i-n,k}a_{k,j}^{\prime}+\sum_{k=m+1}^{m+m^{\prime}}%
d_{i-n,k-m}c_{k-m,j}^{\prime}, & \text{if }i>n\text{ and }j\leq\ell;\\
\sum_{k=1}^{m}c_{i-n,k}b_{k,j-\ell}^{\prime}+\sum_{k=m+1}^{m+m^{\prime}%
}d_{i-n,k-m}d_{k-m,j-\ell}^{\prime}, & \text{if }i>n\text{ and }j>\ell
\end{cases}
\\
&  =\sum_{k=1}^{m}a_{i,k}a_{k,j}^{\prime}+\sum_{k=m+1}^{m+m^{\prime}}%
b_{i,k-m}c_{k-m,j}^{\prime}\ \ \ \ \ \ \ \ \ \ \left(  \text{since }i\leq
n\text{ and }j\leq\ell\right)  ,
\end{align*}
we obtain
\begin{align*}
&  \sum_{k=1}^{m+m^{\prime}}%
\begin{cases}
a_{i,k}, & \text{if }i\leq n\text{ and }k\leq m;\\
b_{i,k-m}, & \text{if }i\leq n\text{ and }k>m;\\
c_{i-n,k}, & \text{if }i>n\text{ and }k\leq m;\\
d_{i-n,k-m}, & \text{if }i>n\text{ and }k>m
\end{cases}%
\begin{cases}
a_{k,j}^{\prime}, & \text{if }k\leq m\text{ and }j\leq\ell;\\
b_{k,j-\ell}^{\prime}, & \text{if }k\leq m\text{ and }j>\ell;\\
c_{k-m,j}^{\prime}, & \text{if }k>m\text{ and }j\leq\ell;\\
d_{k-m,j-\ell}^{\prime}, & \text{if }k>m\text{ and }j>\ell
\end{cases}
\\
&  =%
\begin{cases}
\sum_{k=1}^{m}a_{i,k}a_{k,j}^{\prime}+\sum_{k=m+1}^{m+m^{\prime}}%
b_{i,k-m}c_{k-m,j}^{\prime}, & \text{if }i\leq n\text{ and }j\leq\ell;\\
\sum_{k=1}^{m}a_{i,k}b_{k,j-\ell}^{\prime}+\sum_{k=m+1}^{m+m^{\prime}%
}b_{i,k-m}d_{k-m,j-\ell}^{\prime}, & \text{if }i\leq n\text{ and }j>\ell;\\
\sum_{k=1}^{m}c_{i-n,k}a_{k,j}^{\prime}+\sum_{k=m+1}^{m+m^{\prime}}%
d_{i-n,k-m}c_{k-m,j}^{\prime}, & \text{if }i>n\text{ and }j\leq\ell;\\
\sum_{k=1}^{m}c_{i-n,k}b_{k,j-\ell}^{\prime}+\sum_{k=m+1}^{m+m^{\prime}%
}d_{i-n,k-m}d_{k-m,j-\ell}^{\prime}, & \text{if }i>n\text{ and }j>\ell
\end{cases}
.
\end{align*}
Thus, (\ref{sol.block2x2.mult.entrywise}) is proven in Subcase 1.1.

Let us now consider Subcase 1.2. In this Subcase, we have $j>\ell$. Now,
comparing%
\begin{align*}
&  \sum_{k=1}^{m+m^{\prime}}%
\begin{cases}
a_{i,k}, & \text{if }i\leq n\text{ and }k\leq m;\\
b_{i,k-m}, & \text{if }i\leq n\text{ and }k>m;\\
c_{i-n,k}, & \text{if }i>n\text{ and }k\leq m;\\
d_{i-n,k-m}, & \text{if }i>n\text{ and }k>m
\end{cases}%
\begin{cases}
a_{k,j}^{\prime}, & \text{if }k\leq m\text{ and }j\leq\ell;\\
b_{k,j-\ell}^{\prime}, & \text{if }k\leq m\text{ and }j>\ell;\\
c_{k-m,j}^{\prime}, & \text{if }k>m\text{ and }j\leq\ell;\\
d_{k-m,j-\ell}^{\prime}, & \text{if }k>m\text{ and }j>\ell
\end{cases}
\\
&  =\sum_{k=1}^{m}\underbrace{%
\begin{cases}
a_{i,k}, & \text{if }i\leq n\text{ and }k\leq m;\\
b_{i,k-m}, & \text{if }i\leq n\text{ and }k>m;\\
c_{i-n,k}, & \text{if }i>n\text{ and }k\leq m;\\
d_{i-n,k-m}, & \text{if }i>n\text{ and }k>m
\end{cases}
}_{\substack{=a_{i,k}\\\text{(since }i\leq n\text{ and }k\leq m\text{)}%
}}\underbrace{%
\begin{cases}
a_{k,j}^{\prime}, & \text{if }k\leq m\text{ and }j\leq\ell;\\
b_{k,j-\ell}^{\prime}, & \text{if }k\leq m\text{ and }j>\ell;\\
c_{k-m,j}^{\prime}, & \text{if }k>m\text{ and }j\leq\ell;\\
d_{k-m,j-\ell}^{\prime}, & \text{if }k>m\text{ and }j>\ell
\end{cases}
}_{\substack{=b_{k,j-\ell}^{\prime}\\\text{(since }k\leq m\text{ and }%
j>\ell\text{)}}}\\
&  \ \ \ \ \ \ \ \ \ \ +\sum_{k=m+1}^{m+m^{\prime}}\underbrace{%
\begin{cases}
a_{i,k}, & \text{if }i\leq n\text{ and }k\leq m;\\
b_{i,k-m}, & \text{if }i\leq n\text{ and }k>m;\\
c_{i-n,k}, & \text{if }i>n\text{ and }k\leq m;\\
d_{i-n,k-m}, & \text{if }i>n\text{ and }k>m
\end{cases}
}_{\substack{=b_{i,k-m}\\\text{(since }i\leq n\text{ and }k>m\text{)}%
}}\underbrace{%
\begin{cases}
a_{k,j}^{\prime}, & \text{if }k\leq m\text{ and }j\leq\ell;\\
b_{k,j-\ell}^{\prime}, & \text{if }k\leq m\text{ and }j>\ell;\\
c_{k-m,j}^{\prime}, & \text{if }k>m\text{ and }j\leq\ell;\\
d_{k-m,j-\ell}^{\prime}, & \text{if }k>m\text{ and }j>\ell
\end{cases}
}_{\substack{=d_{k-m,j-\ell}^{\prime}\\\text{(since }k>m\text{ and }%
j>\ell\text{)}}}\\
&  \ \ \ \ \ \ \ \ \ \ \left(  \text{since }0\leq m\leq m+m^{\prime}\right) \\
&  =\sum_{k=1}^{m}a_{i,k}b_{k,j-\ell}^{\prime}+\sum_{k=m+1}^{m+m^{\prime}%
}b_{i,k-m}d_{k-m,j-\ell}^{\prime}%
\end{align*}
with%
\begin{align*}
&
\begin{cases}
\sum_{k=1}^{m}a_{i,k}a_{k,j}^{\prime}+\sum_{k=m+1}^{m+m^{\prime}}%
b_{i,k-m}c_{k-m,j}^{\prime}, & \text{if }i\leq n\text{ and }j\leq\ell;\\
\sum_{k=1}^{m}a_{i,k}b_{k,j-\ell}^{\prime}+\sum_{k=m+1}^{m+m^{\prime}%
}b_{i,k-m}d_{k-m,j-\ell}^{\prime}, & \text{if }i\leq n\text{ and }j>\ell;\\
\sum_{k=1}^{m}c_{i-n,k}a_{k,j}^{\prime}+\sum_{k=m+1}^{m+m^{\prime}}%
d_{i-n,k-m}c_{k-m,j}^{\prime}, & \text{if }i>n\text{ and }j\leq\ell;\\
\sum_{k=1}^{m}c_{i-n,k}b_{k,j-\ell}^{\prime}+\sum_{k=m+1}^{m+m^{\prime}%
}d_{i-n,k-m}d_{k-m,j-\ell}^{\prime}, & \text{if }i>n\text{ and }j>\ell
\end{cases}
\\
&  =\sum_{k=1}^{m}a_{i,k}b_{k,j-\ell}^{\prime}+\sum_{k=m+1}^{m+m^{\prime}%
}b_{i,k-m}d_{k-m,j-\ell}^{\prime}\ \ \ \ \ \ \ \ \ \ \left(  \text{since
}i\leq n\text{ and }j>\ell\right)  ,
\end{align*}
we obtain
\begin{align*}
&  \sum_{k=1}^{m+m^{\prime}}%
\begin{cases}
a_{i,k}, & \text{if }i\leq n\text{ and }k\leq m;\\
b_{i,k-m}, & \text{if }i\leq n\text{ and }k>m;\\
c_{i-n,k}, & \text{if }i>n\text{ and }k\leq m;\\
d_{i-n,k-m}, & \text{if }i>n\text{ and }k>m
\end{cases}%
\begin{cases}
a_{k,j}^{\prime}, & \text{if }k\leq m\text{ and }j\leq\ell;\\
b_{k,j-\ell}^{\prime}, & \text{if }k\leq m\text{ and }j>\ell;\\
c_{k-m,j}^{\prime}, & \text{if }k>m\text{ and }j\leq\ell;\\
d_{k-m,j-\ell}^{\prime}, & \text{if }k>m\text{ and }j>\ell
\end{cases}
\\
&  =%
\begin{cases}
\sum_{k=1}^{m}a_{i,k}a_{k,j}^{\prime}+\sum_{k=m+1}^{m+m^{\prime}}%
b_{i,k-m}c_{k-m,j}^{\prime}, & \text{if }i\leq n\text{ and }j\leq\ell;\\
\sum_{k=1}^{m}a_{i,k}b_{k,j-\ell}^{\prime}+\sum_{k=m+1}^{m+m^{\prime}%
}b_{i,k-m}d_{k-m,j-\ell}^{\prime}, & \text{if }i\leq n\text{ and }j>\ell;\\
\sum_{k=1}^{m}c_{i-n,k}a_{k,j}^{\prime}+\sum_{k=m+1}^{m+m^{\prime}}%
d_{i-n,k-m}c_{k-m,j}^{\prime}, & \text{if }i>n\text{ and }j\leq\ell;\\
\sum_{k=1}^{m}c_{i-n,k}b_{k,j-\ell}^{\prime}+\sum_{k=m+1}^{m+m^{\prime}%
}d_{i-n,k-m}d_{k-m,j-\ell}^{\prime}, & \text{if }i>n\text{ and }j>\ell
\end{cases}
.
\end{align*}
Thus, (\ref{sol.block2x2.mult.entrywise}) is proven in Subcase 1.2.

We have thus proven (\ref{sol.block2x2.mult.entrywise}) in each of the two
Subcases 1.1 and 1.2. Since these two Subcases cover the whole Case 1, this
shows that (\ref{sol.block2x2.mult.entrywise}) is proven in Case 1.

Let us now consider Case 2. In this case, we have $i>n$. Now, we must be in
one of the following two subcases:

\textit{Subcase 2.1:} We have $j\leq\ell$.

\textit{Subcase 2.2:} We have $j>\ell$.

Let us first consider Subcase 2.1. In this Subcase, we have $j\leq\ell$. Now,
comparing%
\begin{align*}
&  \sum_{k=1}^{m+m^{\prime}}%
\begin{cases}
a_{i,k}, & \text{if }i\leq n\text{ and }k\leq m;\\
b_{i,k-m}, & \text{if }i\leq n\text{ and }k>m;\\
c_{i-n,k}, & \text{if }i>n\text{ and }k\leq m;\\
d_{i-n,k-m}, & \text{if }i>n\text{ and }k>m
\end{cases}%
\begin{cases}
a_{k,j}^{\prime}, & \text{if }k\leq m\text{ and }j\leq\ell;\\
b_{k,j-\ell}^{\prime}, & \text{if }k\leq m\text{ and }j>\ell;\\
c_{k-m,j}^{\prime}, & \text{if }k>m\text{ and }j\leq\ell;\\
d_{k-m,j-\ell}^{\prime}, & \text{if }k>m\text{ and }j>\ell
\end{cases}
\\
&  =\sum_{k=1}^{m}\underbrace{%
\begin{cases}
a_{i,k}, & \text{if }i\leq n\text{ and }k\leq m;\\
b_{i,k-m}, & \text{if }i\leq n\text{ and }k>m;\\
c_{i-n,k}, & \text{if }i>n\text{ and }k\leq m;\\
d_{i-n,k-m}, & \text{if }i>n\text{ and }k>m
\end{cases}
}_{\substack{=c_{i-n,k}\\\text{(since }i>n\text{ and }k\leq m\text{)}%
}}\underbrace{%
\begin{cases}
a_{k,j}^{\prime}, & \text{if }k\leq m\text{ and }j\leq\ell;\\
b_{k,j-\ell}^{\prime}, & \text{if }k\leq m\text{ and }j>\ell;\\
c_{k-m,j}^{\prime}, & \text{if }k>m\text{ and }j\leq\ell;\\
d_{k-m,j-\ell}^{\prime}, & \text{if }k>m\text{ and }j>\ell
\end{cases}
}_{\substack{=a_{k,j}^{\prime}\\\text{(since }k\leq m\text{ and }j\leq
\ell\text{)}}}\\
&  \ \ \ \ \ \ \ \ \ \ +\sum_{k=m+1}^{m+m^{\prime}}\underbrace{%
\begin{cases}
a_{i,k}, & \text{if }i\leq n\text{ and }k\leq m;\\
b_{i,k-m}, & \text{if }i\leq n\text{ and }k>m;\\
c_{i-n,k}, & \text{if }i>n\text{ and }k\leq m;\\
d_{i-n,k-m}, & \text{if }i>n\text{ and }k>m
\end{cases}
}_{\substack{=d_{i-n,k-m}\\\text{(since }i>n\text{ and }k>m\text{)}%
}}\underbrace{%
\begin{cases}
a_{k,j}^{\prime}, & \text{if }k\leq m\text{ and }j\leq\ell;\\
b_{k,j-\ell}^{\prime}, & \text{if }k\leq m\text{ and }j>\ell;\\
c_{k-m,j}^{\prime}, & \text{if }k>m\text{ and }j\leq\ell;\\
d_{k-m,j-\ell}^{\prime}, & \text{if }k>m\text{ and }j>\ell
\end{cases}
}_{\substack{=c_{k-m,j}^{\prime}\\\text{(since }k>m\text{ and }j\leq
\ell\text{)}}}\\
&  \ \ \ \ \ \ \ \ \ \ \left(  \text{since }0\leq m\leq m+m^{\prime}\right) \\
&  =\sum_{k=1}^{m}c_{i-n,k}a_{k,j}^{\prime}+\sum_{k=m+1}^{m+m^{\prime}%
}d_{i-n,k-m}c_{k-m,j}^{\prime}%
\end{align*}
with%
\begin{align*}
&
\begin{cases}
\sum_{k=1}^{m}a_{i,k}a_{k,j}^{\prime}+\sum_{k=m+1}^{m+m^{\prime}}%
b_{i,k-m}c_{k-m,j}^{\prime}, & \text{if }i\leq n\text{ and }j\leq\ell;\\
\sum_{k=1}^{m}a_{i,k}b_{k,j-\ell}^{\prime}+\sum_{k=m+1}^{m+m^{\prime}%
}b_{i,k-m}d_{k-m,j-\ell}^{\prime}, & \text{if }i\leq n\text{ and }j>\ell;\\
\sum_{k=1}^{m}c_{i-n,k}a_{k,j}^{\prime}+\sum_{k=m+1}^{m+m^{\prime}}%
d_{i-n,k-m}c_{k-m,j}^{\prime}, & \text{if }i>n\text{ and }j\leq\ell;\\
\sum_{k=1}^{m}c_{i-n,k}b_{k,j-\ell}^{\prime}+\sum_{k=m+1}^{m+m^{\prime}%
}d_{i-n,k-m}d_{k-m,j-\ell}^{\prime}, & \text{if }i>n\text{ and }j>\ell
\end{cases}
\\
&  =\sum_{k=1}^{m}c_{i-n,k}a_{k,j}^{\prime}+\sum_{k=m+1}^{m+m^{\prime}%
}d_{i-n,k-m}c_{k-m,j}^{\prime}\ \ \ \ \ \ \ \ \ \ \left(  \text{since
}i>n\text{ and }j\leq\ell\right)  ,
\end{align*}
we obtain
\begin{align*}
&  \sum_{k=1}^{m+m^{\prime}}%
\begin{cases}
a_{i,k}, & \text{if }i\leq n\text{ and }k\leq m;\\
b_{i,k-m}, & \text{if }i\leq n\text{ and }k>m;\\
c_{i-n,k}, & \text{if }i>n\text{ and }k\leq m;\\
d_{i-n,k-m}, & \text{if }i>n\text{ and }k>m
\end{cases}%
\begin{cases}
a_{k,j}^{\prime}, & \text{if }k\leq m\text{ and }j\leq\ell;\\
b_{k,j-\ell}^{\prime}, & \text{if }k\leq m\text{ and }j>\ell;\\
c_{k-m,j}^{\prime}, & \text{if }k>m\text{ and }j\leq\ell;\\
d_{k-m,j-\ell}^{\prime}, & \text{if }k>m\text{ and }j>\ell
\end{cases}
\\
&  =%
\begin{cases}
\sum_{k=1}^{m}a_{i,k}a_{k,j}^{\prime}+\sum_{k=m+1}^{m+m^{\prime}}%
b_{i,k-m}c_{k-m,j}^{\prime}, & \text{if }i\leq n\text{ and }j\leq\ell;\\
\sum_{k=1}^{m}a_{i,k}b_{k,j-\ell}^{\prime}+\sum_{k=m+1}^{m+m^{\prime}%
}b_{i,k-m}d_{k-m,j-\ell}^{\prime}, & \text{if }i\leq n\text{ and }j>\ell;\\
\sum_{k=1}^{m}c_{i-n,k}a_{k,j}^{\prime}+\sum_{k=m+1}^{m+m^{\prime}}%
d_{i-n,k-m}c_{k-m,j}^{\prime}, & \text{if }i>n\text{ and }j\leq\ell;\\
\sum_{k=1}^{m}c_{i-n,k}b_{k,j-\ell}^{\prime}+\sum_{k=m+1}^{m+m^{\prime}%
}d_{i-n,k-m}d_{k-m,j-\ell}^{\prime}, & \text{if }i>n\text{ and }j>\ell
\end{cases}
.
\end{align*}
Thus, (\ref{sol.block2x2.mult.entrywise}) is proven in Subcase 2.1.

Let us now consider Subcase 2.2. In this Subcase, we have $j>\ell$. Now,
comparing%
\begin{align*}
&  \sum_{k=1}^{m+m^{\prime}}%
\begin{cases}
a_{i,k}, & \text{if }i\leq n\text{ and }k\leq m;\\
b_{i,k-m}, & \text{if }i\leq n\text{ and }k>m;\\
c_{i-n,k}, & \text{if }i>n\text{ and }k\leq m;\\
d_{i-n,k-m}, & \text{if }i>n\text{ and }k>m
\end{cases}%
\begin{cases}
a_{k,j}^{\prime}, & \text{if }k\leq m\text{ and }j\leq\ell;\\
b_{k,j-\ell}^{\prime}, & \text{if }k\leq m\text{ and }j>\ell;\\
c_{k-m,j}^{\prime}, & \text{if }k>m\text{ and }j\leq\ell;\\
d_{k-m,j-\ell}^{\prime}, & \text{if }k>m\text{ and }j>\ell
\end{cases}
\\
&  =\sum_{k=1}^{m}\underbrace{%
\begin{cases}
a_{i,k}, & \text{if }i\leq n\text{ and }k\leq m;\\
b_{i,k-m}, & \text{if }i\leq n\text{ and }k>m;\\
c_{i-n,k}, & \text{if }i>n\text{ and }k\leq m;\\
d_{i-n,k-m}, & \text{if }i>n\text{ and }k>m
\end{cases}
}_{\substack{=c_{i-n,k}\\\text{(since }i>n\text{ and }k\leq m\text{)}%
}}\underbrace{%
\begin{cases}
a_{k,j}^{\prime}, & \text{if }k\leq m\text{ and }j\leq\ell;\\
b_{k,j-\ell}^{\prime}, & \text{if }k\leq m\text{ and }j>\ell;\\
c_{k-m,j}^{\prime}, & \text{if }k>m\text{ and }j\leq\ell;\\
d_{k-m,j-\ell}^{\prime}, & \text{if }k>m\text{ and }j>\ell
\end{cases}
}_{\substack{=b_{k,j-\ell}^{\prime}\\\text{(since }k\leq m\text{ and }%
j>\ell\text{)}}}\\
&  \ \ \ \ \ \ \ \ \ \ +\sum_{k=m+1}^{m+m^{\prime}}\underbrace{%
\begin{cases}
a_{i,k}, & \text{if }i\leq n\text{ and }k\leq m;\\
b_{i,k-m}, & \text{if }i\leq n\text{ and }k>m;\\
c_{i-n,k}, & \text{if }i>n\text{ and }k\leq m;\\
d_{i-n,k-m}, & \text{if }i>n\text{ and }k>m
\end{cases}
}_{\substack{=d_{i-n,k-m}\\\text{(since }i>n\text{ and }k>m\text{)}%
}}\underbrace{%
\begin{cases}
a_{k,j}^{\prime}, & \text{if }k\leq m\text{ and }j\leq\ell;\\
b_{k,j-\ell}^{\prime}, & \text{if }k\leq m\text{ and }j>\ell;\\
c_{k-m,j}^{\prime}, & \text{if }k>m\text{ and }j\leq\ell;\\
d_{k-m,j-\ell}^{\prime}, & \text{if }k>m\text{ and }j>\ell
\end{cases}
}_{\substack{=d_{k-m,j-\ell}^{\prime}\\\text{(since }k>m\text{ and }%
j>\ell\text{)}}}\\
&  \ \ \ \ \ \ \ \ \ \ \left(  \text{since }0\leq m\leq m+m^{\prime}\right) \\
&  =\sum_{k=1}^{m}c_{i-n,k}b_{k,j-\ell}^{\prime}+\sum_{k=m+1}^{m+m^{\prime}%
}d_{i-n,k-m}d_{k-m,j-\ell}^{\prime}%
\end{align*}
with%
\begin{align*}
&
\begin{cases}
\sum_{k=1}^{m}a_{i,k}a_{k,j}^{\prime}+\sum_{k=m+1}^{m+m^{\prime}}%
b_{i,k-m}c_{k-m,j}^{\prime}, & \text{if }i\leq n\text{ and }j\leq\ell;\\
\sum_{k=1}^{m}a_{i,k}b_{k,j-\ell}^{\prime}+\sum_{k=m+1}^{m+m^{\prime}%
}b_{i,k-m}d_{k-m,j-\ell}^{\prime}, & \text{if }i\leq n\text{ and }j>\ell;\\
\sum_{k=1}^{m}c_{i-n,k}a_{k,j}^{\prime}+\sum_{k=m+1}^{m+m^{\prime}}%
d_{i-n,k-m}c_{k-m,j}^{\prime}, & \text{if }i>n\text{ and }j\leq\ell;\\
\sum_{k=1}^{m}c_{i-n,k}b_{k,j-\ell}^{\prime}+\sum_{k=m+1}^{m+m^{\prime}%
}d_{i-n,k-m}d_{k-m,j-\ell}^{\prime}, & \text{if }i>n\text{ and }j>\ell
\end{cases}
\\
&  =\sum_{k=1}^{m}c_{i-n,k}b_{k,j-\ell}^{\prime}+\sum_{k=m+1}^{m+m^{\prime}%
}d_{i-n,k-m}d_{k-m,j-\ell}^{\prime}\ \ \ \ \ \ \ \ \ \ \left(  \text{since
}i>n\text{ and }j>\ell\right)  ,
\end{align*}
we obtain
\begin{align*}
&  \sum_{k=1}^{m+m^{\prime}}%
\begin{cases}
a_{i,k}, & \text{if }i\leq n\text{ and }k\leq m;\\
b_{i,k-m}, & \text{if }i\leq n\text{ and }k>m;\\
c_{i-n,k}, & \text{if }i>n\text{ and }k\leq m;\\
d_{i-n,k-m}, & \text{if }i>n\text{ and }k>m
\end{cases}%
\begin{cases}
a_{k,j}^{\prime}, & \text{if }k\leq m\text{ and }j\leq\ell;\\
b_{k,j-\ell}^{\prime}, & \text{if }k\leq m\text{ and }j>\ell;\\
c_{k-m,j}^{\prime}, & \text{if }k>m\text{ and }j\leq\ell;\\
d_{k-m,j-\ell}^{\prime}, & \text{if }k>m\text{ and }j>\ell
\end{cases}
\\
&  =%
\begin{cases}
\sum_{k=1}^{m}a_{i,k}a_{k,j}^{\prime}+\sum_{k=m+1}^{m+m^{\prime}}%
b_{i,k-m}c_{k-m,j}^{\prime}, & \text{if }i\leq n\text{ and }j\leq\ell;\\
\sum_{k=1}^{m}a_{i,k}b_{k,j-\ell}^{\prime}+\sum_{k=m+1}^{m+m^{\prime}%
}b_{i,k-m}d_{k-m,j-\ell}^{\prime}, & \text{if }i\leq n\text{ and }j>\ell;\\
\sum_{k=1}^{m}c_{i-n,k}a_{k,j}^{\prime}+\sum_{k=m+1}^{m+m^{\prime}}%
d_{i-n,k-m}c_{k-m,j}^{\prime}, & \text{if }i>n\text{ and }j\leq\ell;\\
\sum_{k=1}^{m}c_{i-n,k}b_{k,j-\ell}^{\prime}+\sum_{k=m+1}^{m+m^{\prime}%
}d_{i-n,k-m}d_{k-m,j-\ell}^{\prime}, & \text{if }i>n\text{ and }j>\ell
\end{cases}
.
\end{align*}
Thus, (\ref{sol.block2x2.mult.entrywise}) is proven in Subcase 2.2.

We have thus proven (\ref{sol.block2x2.mult.entrywise}) in each of the two
Subcases 2.1 and 2.2. Since these two Subcases cover the whole Case 2, this
shows that (\ref{sol.block2x2.mult.entrywise}) is proven in Case 2.

We have thus proven (\ref{sol.block2x2.mult.entrywise}) in each of the two
Cases 1 and 2. Thus, (\ref{sol.block2x2.mult.entrywise}) always holds. This
completes our proof of (\ref{sol.block2x2.mult.entrywise}).]

Now, (\ref{sol.block2x2.mult.3}) becomes%
\begin{align*}
&  \left(
\begin{array}
[c]{cc}%
A & B\\
C & D
\end{array}
\right)  \left(
\begin{array}
[c]{cc}%
A^{\prime} & B^{\prime}\\
C^{\prime} & D^{\prime}%
\end{array}
\right) \\
&  =\left(  \underbrace{\sum_{k=1}^{m+m^{\prime}}%
\begin{cases}
a_{i,k}, & \text{if }i\leq n\text{ and }k\leq m;\\
b_{i,k-m}, & \text{if }i\leq n\text{ and }k>m;\\
c_{i-n,k}, & \text{if }i>n\text{ and }k\leq m;\\
d_{i-n,k-m}, & \text{if }i>n\text{ and }k>m
\end{cases}%
\begin{cases}
a_{k,j}^{\prime}, & \text{if }k\leq m\text{ and }j\leq\ell;\\
b_{k,j-\ell}^{\prime}, & \text{if }k\leq m\text{ and }j>\ell;\\
c_{k-m,j}^{\prime}, & \text{if }k>m\text{ and }j\leq\ell;\\
d_{k-m,j-\ell}^{\prime}, & \text{if }k>m\text{ and }j>\ell
\end{cases}
}_{\substack{=%
\begin{cases}
\sum_{k=1}^{m}a_{i,k}a_{k,j}^{\prime}+\sum_{k=m+1}^{m+m^{\prime}}%
b_{i,k-m}c_{k-m,j}^{\prime}, & \text{if }i\leq n\text{ and }j\leq\ell;\\
\sum_{k=1}^{m}a_{i,k}b_{k,j-\ell}^{\prime}+\sum_{k=m+1}^{m+m^{\prime}%
}b_{i,k-m}d_{k-m,j-\ell}^{\prime}, & \text{if }i\leq n\text{ and }j>\ell;\\
\sum_{k=1}^{m}c_{i-n,k}a_{k,j}^{\prime}+\sum_{k=m+1}^{m+m^{\prime}}%
d_{i-n,k-m}c_{k-m,j}^{\prime}, & \text{if }i>n\text{ and }j\leq\ell;\\
\sum_{k=1}^{m}c_{i-n,k}b_{k,j-\ell}^{\prime}+\sum_{k=m+1}^{m+m^{\prime}%
}d_{i-n,k-m}d_{k-m,j-\ell}^{\prime}, & \text{if }i>n\text{ and }j>\ell
\end{cases}
\\\text{(by (\ref{sol.block2x2.mult.entrywise}))}}}\right)  _{1\leq i\leq
n+n^{\prime},\ 1\leq j\leq\ell+\ell^{\prime}}\\
&  =\left(
\begin{cases}
\sum_{k=1}^{m}a_{i,k}a_{k,j}^{\prime}+\sum_{k=m+1}^{m+m^{\prime}}%
b_{i,k-m}c_{k-m,j}^{\prime}, & \text{if }i\leq n\text{ and }j\leq\ell;\\
\sum_{k=1}^{m}a_{i,k}b_{k,j-\ell}^{\prime}+\sum_{k=m+1}^{m+m^{\prime}%
}b_{i,k-m}d_{k-m,j-\ell}^{\prime}, & \text{if }i\leq n\text{ and }j>\ell;\\
\sum_{k=1}^{m}c_{i-n,k}a_{k,j}^{\prime}+\sum_{k=m+1}^{m+m^{\prime}}%
d_{i-n,k-m}c_{k-m,j}^{\prime}, & \text{if }i>n\text{ and }j\leq\ell;\\
\sum_{k=1}^{m}c_{i-n,k}b_{k,j-\ell}^{\prime}+\sum_{k=m+1}^{m+m^{\prime}%
}d_{i-n,k-m}d_{k-m,j-\ell}^{\prime}, & \text{if }i>n\text{ and }j>\ell
\end{cases}
\right)  _{1\leq i\leq n+n^{\prime},\ 1\leq j\leq\ell+\ell^{\prime}}\\
&  =\left(
\begin{array}
[c]{cc}%
AA^{\prime}+BC^{\prime} & AB^{\prime}+BD^{\prime}\\
CA^{\prime}+DC^{\prime} & CB^{\prime}+DD^{\prime}%
\end{array}
\right)  \ \ \ \ \ \ \ \ \ \ \left(  \text{by (\ref{sol.block2x2.mult.9}%
)}\right)  .
\end{align*}
This solves Exercise \ref{exe.block2x2.mult}.
\end{proof}
\end{verlong}

\subsection{Solution to Exercise \ref{exe.block2x2.tridet}}

\begin{vershort}
\begin{proof}
[Solution to Exercise \ref{exe.block2x2.tridet}.]We shall prove Exercise
\ref{exe.block2x2.tridet} by induction over $m$:

\textit{Induction base:} Let $n\in\mathbb{N}$. Let $A$ be an $n\times
n$-matrix. Let $B$ be an $n\times0$-matrix\footnote{Of course, there is only
one such $n\times0$-matrix (namely, the empty matrix).}. Let $D$ be a
$0\times0$-matrix\footnote{Of course, there is only one such $0\times0$-matrix
(namely, the empty matrix).}. Then, all three matrices $B$, $0_{0\times n}$
and $D$ are empty (in the sense that each of them has either $0$ rows or $0$
columns or both), and thus we have $\left(
\begin{array}
[c]{cc}%
A & B\\
0_{0\times n} & D
\end{array}
\right)  =A$. Hence, $\det\left(
\begin{array}
[c]{cc}%
A & B\\
0_{0\times n} & D
\end{array}
\right)  =\det A$. Combined with $\det D=1$ (since $D$ is a $0\times
0$-matrix), this yields $\det\left(
\begin{array}
[c]{cc}%
A & B\\
0_{0\times n} & D
\end{array}
\right)  =\det A\cdot\det D$.

Now, let us forget that we fixed $n$, $A$, $B$ and $D$. We thus have proven
that every $n\in\mathbb{N}$, every $n\times n$-matrix $A$, every $n\times
0$-matrix $B$ and every $0\times0$-matrix $D$ satisfy $\det\left(
\begin{array}
[c]{cc}%
A & B\\
0_{0\times n} & D
\end{array}
\right)  =\det A\cdot\det D$. In other words, Exercise
\ref{exe.block2x2.tridet} holds for $m=0$. This completes the induction base.

\textit{Induction step:} Let $M\in\mathbb{N}$ be positive. Assume that
Exercise \ref{exe.block2x2.tridet} holds for $m=M-1$. We need to prove that
Exercise \ref{exe.block2x2.tridet} holds for $m=M$.

Let $n\in\mathbb{N}$. Let $A$ be an $n\times n$-matrix. Let $B$ be an $n\times
M$-matrix. Let $D$ be an $M\times M$-matrix.

Write the $M\times M$-matrix $D$ in the form $D=\left(  d_{i,j}\right)
_{1\leq i\leq M,\ 1\leq j\leq M}$. Hence, Theorem \ref{thm.laplace.gen}
\textbf{(a)} (applied to $M$, $D$, $d_{i,j}$ and $M$ instead of $n$, $A$,
$a_{i,j}$ and $p$) shows that%
\begin{equation}
\det D=\sum_{q=1}^{M}\left(  -1\right)  ^{M+q}d_{M,q}\det\left(  D_{\sim
M,\sim q}\right)  \label{sol.block2x2.tridet.short.indstep.detD}%
\end{equation}
(since $M\in\left\{  1,2,\ldots,M\right\}  $ (because $M>0$)).

Write the $\left(  n+M\right)  \times\left(  n+M\right)  $-matrix $\left(
\begin{array}
[c]{cc}%
A & B\\
0_{M\times n} & D
\end{array}
\right)  $ in the form $\left(
\begin{array}
[c]{cc}%
A & B\\
0_{M\times n} & D
\end{array}
\right)  =\left(  u_{i,j}\right)  _{1\leq u\leq n+M,\ 1\leq v\leq n+M}$. Thus,%
\begin{equation}
u_{n+M,q}=0\ \ \ \ \ \ \ \ \ \ \text{for every }q\in\left\{  1,2,\ldots
,n\right\}  , \label{sol.block2x2.tridet.short.indstep.u=0}%
\end{equation}
and%
\begin{equation}
u_{n+M,n+q}=d_{M,q}\ \ \ \ \ \ \ \ \ \ \text{for every }q\in\left\{
1,2,\ldots,M\right\}  . \label{sol.block2x2.tridet.short.indstep.u=d}%
\end{equation}

Furthermore, for every $q\in\left\{  1,2,\ldots,M\right\}  $, we have%
\begin{equation}
\left(
\begin{array}
[c]{cc}%
A & B\\
0_{M\times n} & D
\end{array}
\right)  _{\sim\left(  n+M\right)  ,\sim\left(  n+q\right)  }=\left(
\begin{array}
[c]{cc}%
A & B_{q}^{\prime}\\
0_{\left(  M-1\right)  \times n} & D_{\sim M,\sim q}%
\end{array}
\right)  , \label{sol.block2x2.tridet.short.indstep.1}%
\end{equation}
where $B_{q}^{\prime}$ is the result of crossing out the $q$-th column in the
matrix $B$. (Draw the matrices and cross out the appropriate rows and columns
to see why this is true.)

But we have assumed that Exercise \ref{exe.block2x2.tridet} holds for $m=M-1$.
Hence, for every $q\in\left\{  1,2,\ldots,M\right\}  $, we can apply Exercise
\ref{exe.block2x2.tridet} to $M-1$, $B_{q}^{\prime}$ and $D_{\sim M,\sim q}$
instead of $m$, $B$ and $D$. As a result, for every $q\in\left\{
1,2,\ldots,M\right\}  $, we obtain%
\[
\det\left(
\begin{array}
[c]{cc}%
A & B_{q}^{\prime}\\
0_{\left(  M-1\right)  \times n} & D_{\sim M,\sim q}%
\end{array}
\right)  =\det A\cdot\det\left(  D_{\sim M,\sim q}\right)  .
\]
Now, taking determinants in (\ref{sol.block2x2.tridet.short.indstep.1}), we
obtain%
\begin{align}
\det\left(  \left(
\begin{array}
[c]{cc}%
A & B\\
0_{M\times n} & D
\end{array}
\right)  _{\sim\left(  n+M\right)  ,\sim\left(  n+q\right)  }\right)   &
=\det\left(
\begin{array}
[c]{cc}%
A & B_{q}^{\prime}\\
0_{\left(  M-1\right)  \times n} & D_{\sim M,\sim q}%
\end{array}
\right) \nonumber\\
&  =\det A\cdot\det\left(  D_{\sim M,\sim q}\right)  .
\label{sol.block2x2.tridet.short.indstep.2}%
\end{align}

But $n+M>0$ and thus $n+M\in\left\{  1,2,\ldots,n+M\right\}  $. Thus, Theorem
\ref{thm.laplace.gen} \textbf{(a)} (applied to $n+M$, $\left(
\begin{array}
[c]{cc}%
A & B\\
0_{M\times n} & D
\end{array}
\right)  $, $u_{i,j}$ and $n+M$ instead of $n$, $A$, $a_{i,j}$ and $p$) shows
that%
\begin{align*}
&  \det\left(
\begin{array}
[c]{cc}%
A & B\\
0_{M\times n} & D
\end{array}
\right) \\
&  =\sum_{q=1}^{n+M}\left(  -1\right)  ^{\left(  n+M\right)  +q}u_{n+M,q}%
\det\left(  \left(
\begin{array}
[c]{cc}%
A & B\\
0_{M\times n} & D
\end{array}
\right)  _{\sim\left(  n+M\right)  ,\sim q}\right) \\
&  =\sum_{q=1}^{n}\left(  -1\right)  ^{\left(  n+M\right)  +q}%
\underbrace{u_{n+M,q}}_{\substack{=0\\\text{(by
(\ref{sol.block2x2.tridet.short.indstep.u=0}))}}}\det\left(  \left(
\begin{array}
[c]{cc}%
A & B\\
0_{M\times n} & D
\end{array}
\right)  _{\sim\left(  n+M\right)  ,\sim q}\right) \\
&  \ \ \ \ \ \ \ \ \ \ +\sum_{q=n+1}^{n+M}\left(  -1\right)  ^{\left(
n+M\right)  +q}u_{n+M,q}\det\left(  \left(
\begin{array}
[c]{cc}%
A & B\\
0_{M\times n} & D
\end{array}
\right)  _{\sim\left(  n+M\right)  ,\sim q}\right) \\
&  \ \ \ \ \ \ \ \ \ \ \left(  \text{since }0\leq n\leq n+M\right) \\
&  =\underbrace{\sum_{q=1}^{n}\left(  -1\right)  ^{\left(  n+M\right)
+q}0\det\left(  \left(
\begin{array}
[c]{cc}%
A & B\\
0_{M\times n} & D
\end{array}
\right)  _{\sim\left(  n+M\right)  ,\sim q}\right)  }_{=0}\\
&  \ \ \ \ \ \ \ \ \ \ +\sum_{q=n+1}^{n+M}\left(  -1\right)  ^{\left(
n+M\right)  +q}u_{n+M,q}\det\left(  \left(
\begin{array}
[c]{cc}%
A & B\\
0_{M\times n} & D
\end{array}
\right)  _{\sim\left(  n+M\right)  ,\sim q}\right) \\
&  =\sum_{q=n+1}^{n+M}\left(  -1\right)  ^{\left(  n+M\right)  +q}%
u_{n+M,q}\det\left(  \left(
\begin{array}
[c]{cc}%
A & B\\
0_{M\times n} & D
\end{array}
\right)  _{\sim\left(  n+M\right)  ,\sim q}\right) \\
&  =\sum_{q=1}^{M}\underbrace{\left(  -1\right)  ^{\left(  n+M\right)
+\left(  n+q\right)  }}_{\substack{=\left(  -1\right)  ^{M+q}\\\text{(since
}\left(  n+M\right)  +\left(  n+q\right)  \\=2n+M+q\equiv
M+q\operatorname{mod}2\text{)}}}\underbrace{u_{n+M,n+q}}_{\substack{=d_{M,q}%
\\\text{(by (\ref{sol.block2x2.tridet.short.indstep.u=d}))}}}\underbrace{\det
\left(  \left(
\begin{array}
[c]{cc}%
A & B\\
0_{M\times n} & D
\end{array}
\right)  _{\sim\left(  n+M\right)  ,\sim\left(  n+q\right)  }\right)
}_{\substack{=\det A\cdot\det\left(  D_{\sim M,\sim q}\right)  \\\text{(by
(\ref{sol.block2x2.tridet.short.indstep.2}))}}}\\
&  \ \ \ \ \ \ \ \ \ \ \left(  \text{here, we have substituted }n+q\text{ for
}q\text{ in the sum}\right) \\
&  =\sum_{q=1}^{M}\left(  -1\right)  ^{M+q}d_{M,q}\det A\cdot\det\left(
D_{\sim M,\sim q}\right)  =\det A\cdot\underbrace{\sum_{q=1}^{M}\left(
-1\right)  ^{M+q}d_{M,q}\det\left(  D_{\sim M,\sim q}\right)  }%
_{\substack{=\det D\\\text{(by (\ref{sol.block2x2.tridet.short.indstep.detD}%
))}}}\\
&  =\det A\cdot\det D.
\end{align*}

Now, let us forget that we fixed $n$, $A$, $B$ and $D$. We thus have shown
that for every $n\in\mathbb{N}$, for every $n\times n$-matrix $A$, for every
$n\times M$-matrix $B$, and for every $M\times M$-matrix $D$, we have
$\det\left(
\begin{array}
[c]{cc}%
A & B\\
0_{M\times n} & D
\end{array}
\right)  =\det A\cdot\det D$. In other words, Exercise
\ref{exe.block2x2.tridet} holds for $m=M$. This completes the induction step.
Hence, Exercise \ref{exe.block2x2.tridet} is solved by induction.
\end{proof}
\end{vershort}

\begin{verlong}
Before we prove Exercise \ref{exe.block2x2.tridet}, we state a lemma (which is
just a straightforward formalization of an obvious fact):

\begin{lemma}
\label{lem.block2x2.tridet.last-row-minor} Let $n\in\mathbb{N}$ and
$m\in\mathbb{N}$. Let $A=\left(  a_{i,j}\right)  _{1\leq i\leq n,\ 1\leq j\leq
m}$ be an $n\times m$-matrix. Let $q\in\left\{  1,2,\ldots,m\right\}  $. For
every $j\in\mathbb{Z}$ and $r\in\mathbb{Z}$, let $\mathbf{d}_{r}\left(
j\right)  $ be the integer $%
\begin{cases}
j, & \text{if }j<r;\\
j+1, & \text{if }j\geq r
\end{cases}
$.

\textbf{(a)} We have $\operatorname*{cols}\nolimits_{1,2,\ldots,\widehat{q}%
,\ldots,m}A=\left(  a_{i,\mathbf{d}_{q}\left(  j\right)  }\right)  _{1\leq
i\leq n,\ 1\leq j\leq m-1}$.

\textbf{(b)} If $n$ is positive, then $A_{\sim n,\sim q}=\left(
a_{i,\mathbf{d}_{q}\left(  j\right)  }\right)  _{1\leq i\leq n-1,\ 1\leq j\leq
m-1}$.
\end{lemma}

\begin{proof}
[Proof of Lemma \ref{lem.block2x2.tridet.last-row-minor}.]We have
$q\in\left\{  1,2,\ldots,m\right\}  $, so that $1\leq q\leq m$, so that $1\leq
m$. Thus, $m-1\in\mathbb{N}$.

Define an $\left(  m-1\right)  $-tuple $\left(  u_{1},u_{2},\ldots
,u_{m-1}\right)  $ by $\left(  u_{1},u_{2},\ldots,u_{m-1}\right)  =\left(
1,2,\ldots,\widehat{q},\ldots,m\right)  $. Thus, for every $j\in\left\{
1,2,\ldots,m-1\right\}  $, we have%
\begin{equation}
u_{j}=%
\begin{cases}
j, & \text{if }j<q;\\
j+1, & \text{if }j\geq q
\end{cases}
=\mathbf{d}_{q}\left(  j\right)
\label{pf.lem.block2x2.tridet.last-row-minor.uj}%
\end{equation}
(since $\mathbf{d}_{q}\left(  j\right)  =%
\begin{cases}
j, & \text{if }j<q;\\
j+1, & \text{if }j\geq q
\end{cases}
$ (by the definition of $\mathbf{d}_{q}\left(  j\right)  $)).

\textbf{(a)} We have%
\begin{align*}
\operatorname*{cols}\nolimits_{1,2,\ldots,\widehat{q},\ldots,m}A  &
=\operatorname*{cols}\nolimits_{u_{1},u_{2},\ldots,u_{m-1}}A\\
&  \ \ \ \ \ \ \ \ \ \ \left(  \text{since }\left(  1,2,\ldots,\widehat{q}%
,\ldots,m\right)  =\left(  u_{1},u_{2},\ldots,u_{m-1}\right)  \right) \\
&  =\left(  a_{i,u_{y}}\right)  _{1\leq i\leq n,\ 1\leq y\leq m-1}%
\end{align*}
(by the definition of $\operatorname*{cols}\nolimits_{u_{1},u_{2}%
,\ldots,u_{m-1}}A$, since $A=\left(  a_{i,j}\right)  _{1\leq i\leq n,\ 1\leq
j\leq m}$). Thus,%
\[
\operatorname*{cols}\nolimits_{1,2,\ldots,\widehat{q},\ldots,m}A=\left(
a_{i,u_{y}}\right)  _{1\leq i\leq n,\ 1\leq y\leq m-1}=\left(  a_{i,u_{j}%
}\right)  _{1\leq i\leq n,\ 1\leq j\leq m-1}%
\]
(here, we renamed the index $\left(  i,y\right)  $ as $\left(  i,j\right)  $).
Hence,%
\[
\operatorname*{cols}\nolimits_{1,2,\ldots,\widehat{q},\ldots,m}A=\left(
\underbrace{a_{i,u_{j}}}_{\substack{=a_{i,\mathbf{d}_{q}\left(  j\right)
}\\\text{(since }u_{j}=\mathbf{d}_{q}\left(  j\right)  \\\text{(by
(\ref{pf.lem.block2x2.tridet.last-row-minor.uj})))}}}\right)  _{1\leq i\leq
n,\ 1\leq j\leq m-1}=\left(  a_{i,\mathbf{d}_{q}\left(  j\right)  }\right)
_{1\leq i\leq n,\ 1\leq j\leq m-1}.
\]
This proves Lemma \ref{lem.block2x2.tridet.last-row-minor} \textbf{(a)}.

\textbf{(b)} Assume that $n$ is positive. The definition of $A_{\sim n,\sim
q}$ yields%
\begin{align*}
A_{\sim n,\sim q}  &  =\operatorname*{sub}\nolimits_{1,2,\ldots,\widehat{n}%
,\ldots,n}^{1,2,\ldots,\widehat{q},\ldots,m}A\\
&  =\operatorname*{sub}\nolimits_{1,2,\ldots,n-1}^{1,2,\ldots,\widehat{q}%
,\ldots,m}A\ \ \ \ \ \ \ \ \ \ \left(  \text{since }\left(  1,2,\ldots
,\widehat{n},\ldots,n\right)  =\left(  1,2,\ldots,n-1\right)  \right) \\
&  =\operatorname*{sub}\nolimits_{1,2,\ldots,n-1}^{u_{1},u_{2},\ldots,u_{m-1}%
}A\ \ \ \ \ \ \ \ \ \ \left(  \text{since }\left(  1,2,\ldots,\widehat{q}%
,\ldots,m\right)  =\left(  u_{1},u_{2},\ldots,u_{m-1}\right)  \right) \\
&  =\left(  a_{x,u_{y}}\right)  _{1\leq x\leq n-1,\ 1\leq y\leq m-1}%
\end{align*}
(by the definition of $\operatorname*{sub}\nolimits_{1,2,\ldots,n-1}%
^{u_{1},u_{2},\ldots,u_{m-1}}A$, since $A=\left(  a_{i,j}\right)  _{1\leq
i\leq n,\ 1\leq j\leq m}$). Thus,%
\[
A_{\sim n,\sim q}=\left(  a_{x,u_{y}}\right)  _{1\leq x\leq n-1,\ 1\leq y\leq
m-1}=\left(  a_{i,u_{j}}\right)  _{1\leq i\leq n-1,\ 1\leq j\leq m-1}%
\]
(here, we renamed the index $\left(  x,y\right)  $ as $\left(  i,j\right)  $).
Hence,%
\[
A_{\sim n,\sim q}=\left(  \underbrace{a_{i,u_{j}}}_{\substack{=a_{i,\mathbf{d}%
_{q}\left(  j\right)  }\\\text{(since }u_{j}=\mathbf{d}_{q}\left(  j\right)
\\\text{(by (\ref{pf.lem.block2x2.tridet.last-row-minor.uj})))}}}\right)
_{1\leq i\leq n-1,\ 1\leq j\leq m-1}=\left(  a_{i,\mathbf{d}_{q}\left(
j\right)  }\right)  _{1\leq i\leq n-1,\ 1\leq j\leq m-1}.
\]
This proves Lemma \ref{lem.block2x2.tridet.last-row-minor} \textbf{(b)}.
\end{proof}

\begin{proof}
[Solution to Exercise \ref{exe.block2x2.tridet}.]We shall prove Exercise
\ref{exe.block2x2.tridet} by induction over $m$:

\textit{Induction base:} Let $n\in\mathbb{N}$. Let $A$ be an $n\times
n$-matrix. Let $B$ be an $n\times0$-matrix\footnote{Of course, there is only
one such $n\times0$-matrix (namely, the empty matrix).}. Let $D$ be a
$0\times0$-matrix\footnote{Of course, there is only one such $0\times0$-matrix
(namely, the empty matrix).}. Then, the matrix $B$ has $0$ columns (since it
is an $n\times0$-matrix), whereas the matrix $0_{0\times n}$ has $0$ rows
(since it is a $0\times n$-matrix), and the matrix $D$ has $0$ rows (since it
is a $0\times0$-matrix). Thus, each of the three matrices $B$, $D$ and
$0_{0\times n}$ has either $0$ rows or $0$ columns. Therefore%
\begin{equation}
\left(
\begin{array}
[c]{cc}%
A & B\\
0_{0\times n} & D
\end{array}
\right)  =A \label{sol.block2x2.tridet.indbase.1}%
\end{equation}
\footnote{Here is a more formal \textit{proof of
(\ref{sol.block2x2.tridet.indbase.1}):} Write the $n\times n$-matrix $A$ in
the form $A=\left(  a_{i,j}\right)  _{1\leq i\leq n,\ 1\leq j\leq n}$.
\par
Write the $n\times0$-matrix $B$ in the form $B=\left(  b_{i,j}\right)  _{1\leq
i\leq n,\ 1\leq j\leq0}$.
\par
We have $0_{0\times n}=\left(  0\right)  _{1\leq i\leq0,\ 1\leq j\leq n}$ (by
the definition of $0_{0\times n}$).
\par
Write the $0\times0$-matrix $D$ in the form $D=\left(  d_{i,j}\right)  _{1\leq
i\leq0,\ 1\leq j\leq0}$.
\par
Now, (\ref{eq.def.block2x2.formal}) (applied to $n$, $0$, $0$, $0_{0\times n}$
and $0$ instead of $n^{\prime}$, $m$, $m^{\prime}$, $C$ and $c_{i,j}$) shows
that
\begin{align*}
\left(
\begin{array}
[c]{cc}%
A & B\\
0_{0\times n} & D
\end{array}
\right)   &  =\left(
\begin{cases}
a_{i,j} & \text{if }i\leq n\text{ and }j\leq n;\\
b_{i,j-n}, & \text{if }i\leq n\text{ and }j>n;\\
0, & \text{if }i>n\text{ and }j\leq n;\\
d_{i-n,j-n}, & \text{if }i>n\text{ and }j>n
\end{cases}
\right)  _{1\leq i\leq n+0,\ 1\leq j\leq n+0}\\
&  =\left(  \underbrace{%
\begin{cases}
a_{i,j} & \text{if }i\leq n\text{ and }j\leq n;\\
b_{i,j-n}, & \text{if }i\leq n\text{ and }j>n;\\
0, & \text{if }i>n\text{ and }j\leq n;\\
d_{i-n,j-n}, & \text{if }i>n\text{ and }j>n
\end{cases}
}_{\substack{=a_{i,j}\\\text{(since }i\leq n\text{ and }j\leq n\text{)}%
}}\right)  _{1\leq i\leq n,\ 1\leq j\leq n}\ \ \ \ \ \ \ \ \ \ \left(
\text{since }n+0=n\right) \\
&  =\left(  a_{i,j}\right)  _{1\leq i\leq n,\ 1\leq j\leq n}=A.
\end{align*}
This proves (\ref{sol.block2x2.tridet.indbase.1}).}. Hence,
\begin{equation}
\det\underbrace{\left(
\begin{array}
[c]{cc}%
A & B\\
0_{0\times n} & D
\end{array}
\right)  }_{=A}=\det A. \label{sol.block2x2.tridet.indbase.2}%
\end{equation}
But $D$ is a $0\times0$-matrix, and thus has determinant $\det D=1$. Hence,
$\det A\cdot\underbrace{\det D}_{=1}=\det A$. Compared with
(\ref{sol.block2x2.tridet.indbase.2}), this yields $\det\left(
\begin{array}
[c]{cc}%
A & B\\
0_{0\times n} & D
\end{array}
\right)  =\det A\cdot\det D$.

Now, let us forget that we fixed $n$, $A$, $B$ and $D$. We thus have proven
that every $n\in\mathbb{N}$, every $n\times n$-matrix $A$, every $n\times
0$-matrix $B$ and every $0\times0$-matrix $D$ satisfy $\det\left(
\begin{array}
[c]{cc}%
A & B\\
0_{0\times n} & D
\end{array}
\right)  =\det A\cdot\det D$. In other words, Exercise
\ref{exe.block2x2.tridet} holds for $m=0$. This completes the induction base.

\textit{Induction step:} Let $M\in\mathbb{N}$ be positive. Assume that
Exercise \ref{exe.block2x2.tridet} holds for $m=M-1$. We need to prove that
Exercise \ref{exe.block2x2.tridet} holds for $m=M$.

Let $n\in\mathbb{N}$. Let $A$ be an $n\times n$-matrix. Let $B$ be an $n\times
M$-matrix. Let $D$ be an $M\times M$-matrix. Let $U$ be the $\left(
n+M\right)  \times\left(  n+M\right)  $-matrix $\left(
\begin{array}
[c]{cc}%
A & B\\
0_{M\times n} & D
\end{array}
\right)  $. Write the matrix $U$ in the form $U=\left(  u_{i,j}\right)
_{1\leq u\leq n+M,\ 1\leq v\leq n+M}$.

Write the $n\times n$-matrix $A$ in the form $A=\left(  a_{i,j}\right)
_{1\leq i\leq n,\ 1\leq j\leq n}$.

Write the $n\times M$-matrix $B$ in the form $B=\left(  b_{i,j}\right)
_{1\leq i\leq n,\ 1\leq j\leq M}$.

We have $0_{M\times n}=\left(  0\right)  _{1\leq i\leq M,\ 1\leq j\leq n}$ (by
the definition of $0_{M\times n}$).

Write the $M\times M$-matrix $D$ in the form $D=\left(  d_{i,j}\right)
_{1\leq i\leq M,\ 1\leq j\leq M}$.

We have $A=\left(  a_{i,j}\right)  _{1\leq i\leq n,\ 1\leq j\leq n}$,
$B=\left(  b_{i,j}\right)  _{1\leq i\leq n,\ 1\leq j\leq M}$, $0_{M\times
n}=\left(  0\right)  _{1\leq i\leq M,\ 1\leq j\leq n}$ and $D=\left(
d_{i,j}\right)  _{1\leq i\leq M,\ 1\leq j\leq M}$. Thus,
(\ref{eq.def.block2x2.formal}) (applied to $n$, $M$, $n$, $M$, $A$, $B$,
$0_{M\times n}$, $D$, $a_{i,j}$, $b_{i,j}$, $0$ and $d_{i,j}$ instead of $n$,
$n^{\prime}$, $m$, $m^{\prime}$, $A$, $B$, $C$, $D$, $a_{i,j}$, $b_{i,j}$,
$c_{i,j}$ and $d_{i,j}$) shows that
\[
\left(
\begin{array}
[c]{cc}%
A & B\\
0_{M\times n} & D
\end{array}
\right)  =\left(
\begin{cases}
a_{i,j} & \text{if }i\leq n\text{ and }j\leq n;\\
b_{i,j-n}, & \text{if }i\leq n\text{ and }j>n;\\
0, & \text{if }i>n\text{ and }j\leq n;\\
d_{i-n,j-n}, & \text{if }i>n\text{ and }j>n
\end{cases}
\right)  _{1\leq i\leq n+M,\ 1\leq j\leq n+M}.
\]
Thus,%
\begin{align}
U  &  =\left(
\begin{array}
[c]{cc}%
A & B\\
0_{M\times n} & D
\end{array}
\right) \nonumber\\
&  =\left(
\begin{cases}
a_{i,j} & \text{if }i\leq n\text{ and }j\leq n;\\
b_{i,j-n}, & \text{if }i\leq n\text{ and }j>n;\\
0, & \text{if }i>n\text{ and }j\leq n;\\
d_{i-n,j-n}, & \text{if }i>n\text{ and }j>n
\end{cases}
\right)  _{1\leq i\leq n+M,\ 1\leq j\leq n+M}.
\label{sol.block2x2.tridet.indstep.U}%
\end{align}

Now, for every $q\in\left\{  1,2,\ldots,M\right\}  $, let us denote the
$n\times\left(  M-1\right)  $-matrix \newline$\operatorname*{cols}%
\nolimits_{1,2,\ldots,\widehat{q},\ldots,M}B$ by $B_{q}^{\prime}$. Thus,%
\[
B_{q}^{\prime}=\operatorname*{cols}\nolimits_{1,2,\ldots,\widehat{q},\ldots
,M}B\ \ \ \ \ \ \ \ \ \ \text{for every }q\in\left\{  1,2,\ldots,M\right\}  .
\]

Write the $\left(  n+M\right)  \times\left(  n+M\right)  $-matrix $U$ in the
form $U=\left(  u_{i,j}\right)  _{1\leq i\leq n+M,\ 1\leq j\leq n+M}$. Thus,%
\[
\left(  u_{i,j}\right)  _{1\leq i\leq n+M,\ 1\leq j\leq n+M}=U=\left(
\begin{cases}
a_{i,j} & \text{if }i\leq n\text{ and }j\leq n;\\
b_{i,j-n}, & \text{if }i\leq n\text{ and }j>n;\\
0, & \text{if }i>n\text{ and }j\leq n;\\
d_{i-n,j-n}, & \text{if }i>n\text{ and }j>n
\end{cases}
\right)  _{1\leq i\leq n+M,\ 1\leq j\leq n+M}%
\]
(by (\ref{sol.block2x2.tridet.indstep.U})). In other words,%
\begin{equation}
u_{i,j}=%
\begin{cases}
a_{i,j} & \text{if }i\leq n\text{ and }j\leq n;\\
b_{i,j-n}, & \text{if }i\leq n\text{ and }j>n;\\
0, & \text{if }i>n\text{ and }j\leq n;\\
d_{i-n,j-n}, & \text{if }i>n\text{ and }j>n
\end{cases}
\label{sol.block2x2.tridet.indstep.minor.pf.1}%
\end{equation}
for all $\left(  i,j\right)  \in\left\{  1,2,\ldots,n+M\right\}  ^{2}$.

Now, it is easy to see that%
\begin{equation}
U_{\sim\left(  n+M\right)  ,\sim\left(  n+q\right)  }=\left(
\begin{array}
[c]{cc}%
A & B_{q}^{\prime}\\
0_{\left(  M-1\right)  \times n} & D_{\sim M,\sim q}%
\end{array}
\right)  \label{sol.block2x2.tridet.indstep.minor}%
\end{equation}
for every $q\in\left\{  1,2,\ldots,M\right\}  $\ \ \ \ \footnote{\textit{Proof
of (\ref{sol.block2x2.tridet.indstep.minor}):} Let $q\in\left\{
1,2,\ldots,M\right\}  $. Notice that $M$ is positive, and thus $n+M$ is
positive. Also, $q\in\left\{  1,2,\ldots,M\right\}  $, so that $n+q\in\left\{
n+1,n+2,\ldots,n+M\right\}  \subseteq\left\{  1,2,\ldots,n+M\right\}  $.
\par
For every $j\in\mathbb{Z}$ and $r\in\mathbb{Z}$, let $\mathbf{d}_{r}\left(
j\right)  $ be the integer $%
\begin{cases}
j, & \text{if }j<r;\\
j+1, & \text{if }j\geq r
\end{cases}
$.
\par
Now, Lemma \ref{lem.block2x2.tridet.last-row-minor} \textbf{(b)} (applied to
$n+M$, $n+M$, $n+q$, $U$ and $u_{i,j}$ instead of $n$, $m$, $q$, $A$ and
$a_{i,j}$) yields%
\begin{equation}
U_{\sim\left(  n+M\right)  ,\sim\left(  n+q\right)  }=\left(  u_{i,\mathbf{d}%
_{n+q}\left(  j\right)  }\right)  _{1\leq i\leq n+M-1,\ 1\leq j\leq n+M-1}
\label{sol.block2x2.tridet.indstep.minor.pf.2}%
\end{equation}
(since $U=\left(  u_{i,j}\right)  _{1\leq i\leq n+M,\ 1\leq j\leq n+M}$).
\par
On the other hand, $D=\left(  d_{i,j}\right)  _{1\leq i\leq M,\ 1\leq j\leq
M}$. Hence, Lemma \ref{lem.block2x2.tridet.last-row-minor} \textbf{(b)}
(applied to $M$, $M$, $D$ and $d_{i,j}$ instead of $n$, $m$, $A$ and $a_{i,j}%
$) yields%
\begin{equation}
D_{\sim M,\sim q}=\left(  d_{i,\mathbf{d}_{q}\left(  j\right)  }\right)
_{1\leq i\leq M-1,\ 1\leq j\leq M-1}.
\label{sol.block2x2.tridet.indstep.minor.pf.3}%
\end{equation}
\par
Also, $B=\left(  b_{i,j}\right)  _{1\leq i\leq n,\ 1\leq j\leq M}$. Hence,
Lemma \ref{lem.block2x2.tridet.last-row-minor} \textbf{(a)} (applied to $n$,
$M$, $B$ and $b_{i,j}$ instead of $n$, $m$, $A$ and $a_{i,j}$) yields%
\[
\operatorname*{cols}\nolimits_{1,2,\ldots,\widehat{q},\ldots,M}B=\left(
b_{i,\mathbf{d}_{q}\left(  j\right)  }\right)  _{1\leq i\leq n,\ 1\leq j\leq
M-1}.
\]
Thus,%
\begin{equation}
B_{q}^{\prime}=\operatorname*{cols}\nolimits_{1,2,\ldots,\widehat{q},\ldots
,M}B=\left(  b_{i,\mathbf{d}_{q}\left(  j\right)  }\right)  _{1\leq i\leq
n,\ 1\leq j\leq M-1}. \label{sol.block2x2.tridet.indstep.minor.pf.4}%
\end{equation}
\par
Now, $A=\left(  a_{i,j}\right)  _{1\leq i\leq n,\ 1\leq j\leq n}$,
$B_{q}^{\prime}=\left(  b_{i,\mathbf{d}_{q}\left(  j\right)  }\right)  _{1\leq
i\leq n,\ 1\leq j\leq M-1}$, $0_{\left(  M-1\right)  \times n}=\left(
0\right)  _{1\leq i\leq M-1,\ 1\leq j\leq n}$ (by the definition of
$0_{\left(  M-1\right)  \times n}$) and $D_{\sim M,\sim q}=\left(
d_{i,\mathbf{d}_{q}\left(  j\right)  }\right)  _{1\leq i\leq M-1,\ 1\leq j\leq
M-1}$. Thus, (\ref{eq.def.block2x2.formal}) (applied to $n$, $M-1$, $n$,
$M-1$, $A$, $B_{q}^{\prime}$, $0_{\left(  M-1\right)  \times n}$, $D_{\sim
M,\sim q}$, $a_{i,j}$, $b_{i,\mathbf{d}_{q}\left(  j\right)  }$, $0$ and
$d_{i,\mathbf{d}_{q}\left(  j\right)  }$ instead of $n$, $n^{\prime}$, $m$,
$m^{\prime}$, $A$, $B$, $C$, $D$, $a_{i,j}$, $b_{i,j}$, $c_{i,j}$ and
$d_{i,j}$) shows that
\begin{equation}
\left(
\begin{array}
[c]{cc}%
A & B_{q}^{\prime}\\
0_{\left(  M-1\right)  \times n} & D_{\sim M,\sim q}%
\end{array}
\right)  =\left(
\begin{cases}
a_{i,j} & \text{if }i\leq n\text{ and }j\leq n;\\
b_{i,\mathbf{d}_{q}\left(  j-n\right)  }, & \text{if }i\leq n\text{ and
}j>n;\\
0, & \text{if }i>n\text{ and }j\leq n;\\
d_{i-n,\mathbf{d}_{q}\left(  j-n\right)  }, & \text{if }i>n\text{ and }j>n
\end{cases}
\right)  _{1\leq i\leq n+M-1,\ 1\leq j\leq n+M-1}.
\label{sol.block2x2.tridet.indstep.minor.pf.6}%
\end{equation}
\par
Now, we shall prove that%
\begin{equation}
u_{i,\mathbf{d}_{n+q}\left(  j\right)  }=%
\begin{cases}
a_{i,j} & \text{if }i\leq n\text{ and }j\leq n;\\
b_{i,\mathbf{d}_{q}\left(  j-n\right)  }, & \text{if }i\leq n\text{ and
}j>n;\\
0, & \text{if }i>n\text{ and }j\leq n;\\
d_{i-n,\mathbf{d}_{q}\left(  j-n\right)  }, & \text{if }i>n\text{ and }j>n
\end{cases}
\label{sol.block2x2.tridet.indstep.minor.pf.8}%
\end{equation}
for every $\left(  i,j\right)  \in\left\{  1,2,\ldots,n+M-1\right\}  ^{2}$.
\par
[\textit{Proof of (\ref{sol.block2x2.tridet.indstep.minor.pf.8}):} Let
$\left(  i,j\right)  \in\left\{  1,2,\ldots,n+M-1\right\}  ^{2}$. Thus,
$i\in\left\{  1,2,\ldots,n+M-1\right\}  $ and $j\in\left\{  1,2,\ldots
,n+M-1\right\}  $.
\par
The definition of $\mathbf{d}_{q}\left(  j-n\right)  $ yields $\mathbf{d}%
_{q}\left(  j-n\right)  =%
\begin{cases}
j-n, & \text{if }j-n<q;\\
j-n+1, & \text{if }j-n\geq q
\end{cases}
$.
\par
The definition of $\mathbf{d}_{q}\left(  j\right)  $ yields $\mathbf{d}%
_{q}\left(  j\right)  =%
\begin{cases}
j, & \text{if }j<q;\\
j+1, & \text{if }j\geq q
\end{cases}
$.
\par
The definition of $\mathbf{d}_{n+q}\left(  j\right)  $ yields%
\begin{align*}
\mathbf{d}_{n+q}\left(  j\right)   &  =%
\begin{cases}
j, & \text{if }j<n+q;\\
j+1, & \text{if }j\geq n+q
\end{cases}
=%
\begin{cases}
j, & \text{if }j-n<q;\\
j+1, & \text{if }j-n\geq q
\end{cases}
\\
&  \ \ \ \ \ \ \ \ \ \ \left(
\begin{array}
[c]{c}%
\text{since the condition }j<n+q\text{ is equivalent to }j-n<q\text{,}\\
\text{and since the condition }j\geq n+q\text{ is equivalent to }j-n\geq q
\end{array}
\right)  .
\end{align*}
Subtracting $n$ from both sides of this equality, we obtain
\begin{align*}
\mathbf{d}_{n+q}\left(  j\right)  -n  &  =%
\begin{cases}
j, & \text{if }j-n<q;\\
j+1, & \text{if }j-n\geq q
\end{cases}
-n=%
\begin{cases}
j-n, & \text{if }j-n<q;\\
j+1-n, & \text{if }j-n\geq q
\end{cases}
\\
&  =%
\begin{cases}
j-n, & \text{if }j-n<q;\\
j-n+1, & \text{if }j-n\geq q
\end{cases}
\ \ \ \ \ \ \ \ \ \ \left(
\begin{array}
[c]{c}%
\text{since }j+1-n=j-n+1\text{ in}\\
\text{the case when }j-n\geq q
\end{array}
\right) \\
&  =\mathbf{d}_{q}\left(  j-n\right)  .
\end{align*}
\par
Notice that%
\begin{align*}
\mathbf{d}_{n+q}\left(  j\right)   &  =%
\begin{cases}
j, & \text{if }j<n+q;\\
j+1, & \text{if }j\geq n+q
\end{cases}
\\
&  \geq%
\begin{cases}
j, & \text{if }j<n+q;\\
j, & \text{if }j\geq n+q
\end{cases}
\ \ \ \ \ \ \ \ \ \ \left(  \text{since }j+1\geq j\text{ in the case when
}j\geq n+q\right) \\
&  =j.
\end{align*}
\par
We must be in one of the following two cases:
\par
\textit{Case 1:} We have $i\leq n$.
\par
\textit{Case 2:} We have $i>n$.
\par
Let us first consider Case 1. In this case, we have $i\leq n$. We are in one
of the following two subcases:
\par
\textit{Subcase 1.1:} We have $j\leq n$.
\par
\textit{Subcase 1.2:} We have $j>n$.
\par
Let us first consider Subcase 1.1. In this Subcase, we have $j\leq n$. Now,
$j\leq n<n+q$ (since $n+\underbrace{q}_{>0}>n$) and $\mathbf{d}_{n+q}\left(
j\right)  =%
\begin{cases}
j, & \text{if }j<n+q;\\
j+1, & \text{if }j\geq n+q
\end{cases}
=j$ (since $j<n+q$). Thus,%
\begin{align*}
u_{i,\mathbf{d}_{n+q}\left(  j\right)  }  &  =u_{i,j}%
\ \ \ \ \ \ \ \ \ \ \left(  \text{since }\mathbf{d}_{n+q}\left(  j\right)
=j\right) \\
&  =%
\begin{cases}
a_{i,j} & \text{if }i\leq n\text{ and }j\leq n;\\
b_{i,j-n}, & \text{if }i\leq n\text{ and }j>n;\\
0, & \text{if }i>n\text{ and }j\leq n;\\
d_{i-n,j-n}, & \text{if }i>n\text{ and }j>n
\end{cases}
\ \ \ \ \ \ \ \ \ \ \left(  \text{by
(\ref{sol.block2x2.tridet.indstep.minor.pf.1})}\right) \\
&  =a_{i,j}%
\end{align*}
(since $i\leq n$ and $j\leq n$). Compared with%
\[%
\begin{cases}
a_{i,j} & \text{if }i\leq n\text{ and }j\leq n;\\
b_{i,\mathbf{d}_{q}\left(  j-n\right)  }, & \text{if }i\leq n\text{ and
}j>n;\\
0, & \text{if }i>n\text{ and }j\leq n;\\
d_{i-n,\mathbf{d}_{q}\left(  j-n\right)  }, & \text{if }i>n\text{ and }j>n
\end{cases}
=a_{i,j}\ \ \ \ \ \ \ \ \ \ \left(  \text{since }i\leq n\text{ and }j\leq
n\right)  ,
\]
this shows that $u_{i,\mathbf{d}_{n+q}\left(  j\right)  }=%
\begin{cases}
a_{i,j} & \text{if }i\leq n\text{ and }j\leq n;\\
b_{i,\mathbf{d}_{q}\left(  j-n\right)  }, & \text{if }i\leq n\text{ and
}j>n;\\
0, & \text{if }i>n\text{ and }j\leq n;\\
d_{i-n,\mathbf{d}_{q}\left(  j-n\right)  }, & \text{if }i>n\text{ and }j>n
\end{cases}
$. Hence, (\ref{sol.block2x2.tridet.indstep.minor.pf.8}) is proven in Subcase
1.1.
\par
Let us now consider Subcase 1.2. In this Subcase, we have $j>n$. Hence,
$\mathbf{d}_{n+q}\left(  j\right)  \geq j>n$. Now,%
\begin{align*}
u_{i,\mathbf{d}_{n+q}\left(  j\right)  }  &  =%
\begin{cases}
a_{i,\mathbf{d}_{n+q}\left(  j\right)  } & \text{if }i\leq n\text{ and
}\mathbf{d}_{n+q}\left(  j\right)  \leq n;\\
b_{i,\mathbf{d}_{n+q}\left(  j\right)  -n}, & \text{if }i\leq n\text{ and
}\mathbf{d}_{n+q}\left(  j\right)  >n;\\
0, & \text{if }i>n\text{ and }\mathbf{d}_{n+q}\left(  j\right)  \leq n;\\
d_{i-n,\mathbf{d}_{n+q}\left(  j\right)  -n}, & \text{if }i>n\text{ and
}\mathbf{d}_{n+q}\left(  j\right)  >n
\end{cases}
\\
&  \ \ \ \ \ \ \ \ \ \ \left(  \text{by
(\ref{sol.block2x2.tridet.indstep.minor.pf.1}), applied to }\left(
i,\mathbf{d}_{n+q}\left(  j\right)  \right)  \text{ instead of }\left(
i,j\right)  \right) \\
&  =b_{i,\mathbf{d}_{n+q}\left(  j\right)  -n}\ \ \ \ \ \ \ \ \ \ \left(
\text{since }i\leq n\text{ and }\mathbf{d}_{n+q}\left(  j\right)  >n\right) \\
&  =b_{i,\mathbf{d}_{q}\left(  j-n\right)  }\ \ \ \ \ \ \ \ \ \ \left(
\text{since }\mathbf{d}_{n+q}\left(  j\right)  -n=\mathbf{d}_{q}\left(
j-n\right)  \right)  .
\end{align*}
Compared with%
\[%
\begin{cases}
a_{i,j} & \text{if }i\leq n\text{ and }j\leq n;\\
b_{i,\mathbf{d}_{q}\left(  j-n\right)  }, & \text{if }i\leq n\text{ and
}j>n;\\
0, & \text{if }i>n\text{ and }j\leq n;\\
d_{i-n,\mathbf{d}_{q}\left(  j-n\right)  }, & \text{if }i>n\text{ and }j>n
\end{cases}
=b_{i,\mathbf{d}_{q}\left(  j-n\right)  }\ \ \ \ \ \ \ \ \ \ \left(
\text{since }i\leq n\text{ and }j>n\right)  ,
\]
this shows that $u_{i,\mathbf{d}_{n+q}\left(  j\right)  }=%
\begin{cases}
a_{i,j} & \text{if }i\leq n\text{ and }j\leq n;\\
b_{i,\mathbf{d}_{q}\left(  j-n\right)  }, & \text{if }i\leq n\text{ and
}j>n;\\
0, & \text{if }i>n\text{ and }j\leq n;\\
d_{i-n,\mathbf{d}_{q}\left(  j-n\right)  }, & \text{if }i>n\text{ and }j>n
\end{cases}
$. Hence, (\ref{sol.block2x2.tridet.indstep.minor.pf.8}) is proven in Subcase
1.2.
\par
We have thus proven (\ref{sol.block2x2.tridet.indstep.minor.pf.8}) in each of
the two Subcases 1.1 and 1.2. Since these two Subcases cover the whole Case 1,
this shows that (\ref{sol.block2x2.tridet.indstep.minor.pf.8}) holds in Case
1.
\par
Let us now consider Case 2. In this case, we have $i>n$. We are in one of the
following two subcases:
\par
\textit{Subcase 2.1:} We have $j\leq n$.
\par
\textit{Subcase 2.2:} We have $j>n$.
\par
Let us first consider Subcase 2.1. In this Subcase, we have $j\leq n$. Now,
$j\leq n<n+q$ (since $n+\underbrace{q}_{>0}>n$) and $\mathbf{d}_{n+q}\left(
j\right)  =%
\begin{cases}
j, & \text{if }j<n+q;\\
j+1, & \text{if }j\geq n+q
\end{cases}
=j$ (since $j<n+q$). Thus,%
\begin{align*}
u_{i,\mathbf{d}_{n+q}\left(  j\right)  }  &  =u_{i,j}%
\ \ \ \ \ \ \ \ \ \ \left(  \text{since }\mathbf{d}_{n+q}\left(  j\right)
=j\right) \\
&  =%
\begin{cases}
a_{i,j} & \text{if }i\leq n\text{ and }j\leq n;\\
b_{i,j-n}, & \text{if }i\leq n\text{ and }j>n;\\
0, & \text{if }i>n\text{ and }j\leq n;\\
d_{i-n,j-n}, & \text{if }i>n\text{ and }j>n
\end{cases}
\ \ \ \ \ \ \ \ \ \ \left(  \text{by
(\ref{sol.block2x2.tridet.indstep.minor.pf.1})}\right) \\
&  =0
\end{align*}
(since $i>n$ and $j\leq n$). Compared with%
\[%
\begin{cases}
a_{i,j} & \text{if }i\leq n\text{ and }j\leq n;\\
b_{i,\mathbf{d}_{q}\left(  j-n\right)  }, & \text{if }i\leq n\text{ and
}j>n;\\
0, & \text{if }i>n\text{ and }j\leq n;\\
d_{i-n,\mathbf{d}_{q}\left(  j-n\right)  }, & \text{if }i>n\text{ and }j>n
\end{cases}
=0\ \ \ \ \ \ \ \ \ \ \left(  \text{since }i>n\text{ and }j\leq n\right)  ,
\]
this shows that $u_{i,\mathbf{d}_{n+q}\left(  j\right)  }=%
\begin{cases}
a_{i,j} & \text{if }i\leq n\text{ and }j\leq n;\\
b_{i,\mathbf{d}_{q}\left(  j-n\right)  }, & \text{if }i\leq n\text{ and
}j>n;\\
0, & \text{if }i>n\text{ and }j\leq n;\\
d_{i-n,\mathbf{d}_{q}\left(  j-n\right)  }, & \text{if }i>n\text{ and }j>n
\end{cases}
$. Hence, (\ref{sol.block2x2.tridet.indstep.minor.pf.8}) is proven in Subcase
2.1.
\par
Let us now consider Subcase 2.2. In this Subcase, we have $j>n$. Hence,
$\mathbf{d}_{n+q}\left(  j\right)  \geq j>n$. Now,%
\begin{align*}
u_{i,\mathbf{d}_{n+q}\left(  j\right)  }  &  =%
\begin{cases}
a_{i,\mathbf{d}_{n+q}\left(  j\right)  } & \text{if }i\leq n\text{ and
}\mathbf{d}_{n+q}\left(  j\right)  \leq n;\\
b_{i,\mathbf{d}_{n+q}\left(  j\right)  -n}, & \text{if }i\leq n\text{ and
}\mathbf{d}_{n+q}\left(  j\right)  >n;\\
0, & \text{if }i>n\text{ and }\mathbf{d}_{n+q}\left(  j\right)  \leq n;\\
d_{i-n,\mathbf{d}_{n+q}\left(  j\right)  -n}, & \text{if }i>n\text{ and
}\mathbf{d}_{n+q}\left(  j\right)  >n
\end{cases}
\\
&  \ \ \ \ \ \ \ \ \ \ \left(  \text{by
(\ref{sol.block2x2.tridet.indstep.minor.pf.1}), applied to }\left(
i,\mathbf{d}_{n+q}\left(  j\right)  \right)  \text{ instead of }\left(
i,j\right)  \right) \\
&  =d_{i-n,\mathbf{d}_{n+q}\left(  j\right)  -n}\ \ \ \ \ \ \ \ \ \ \left(
\text{since }i>n\text{ and }\mathbf{d}_{n+q}\left(  j\right)  >n\right) \\
&  =d_{i-n,\mathbf{d}_{q}\left(  j-n\right)  }\ \ \ \ \ \ \ \ \ \ \left(
\text{since }\mathbf{d}_{n+q}\left(  j\right)  -n=\mathbf{d}_{q}\left(
j-n\right)  \right)  .
\end{align*}
Compared with%
\[%
\begin{cases}
a_{i,j} & \text{if }i\leq n\text{ and }j\leq n;\\
b_{i,\mathbf{d}_{q}\left(  j-n\right)  }, & \text{if }i\leq n\text{ and
}j>n;\\
0, & \text{if }i>n\text{ and }j\leq n;\\
d_{i-n,\mathbf{d}_{q}\left(  j-n\right)  }, & \text{if }i>n\text{ and }j>n
\end{cases}
=d_{i-n,\mathbf{d}_{q}\left(  j-n\right)  }\ \ \ \ \ \ \ \ \ \ \left(
\text{since }i>n\text{ and }j>n\right)  ,
\]
this shows that $u_{i,\mathbf{d}_{n+q}\left(  j\right)  }=%
\begin{cases}
a_{i,j} & \text{if }i\leq n\text{ and }j\leq n;\\
b_{i,\mathbf{d}_{q}\left(  j-n\right)  }, & \text{if }i\leq n\text{ and
}j>n;\\
0, & \text{if }i>n\text{ and }j\leq n;\\
d_{i-n,\mathbf{d}_{q}\left(  j-n\right)  }, & \text{if }i>n\text{ and }j>n
\end{cases}
$. Hence, (\ref{sol.block2x2.tridet.indstep.minor.pf.8}) is proven in Subcase
2.2.
\par
We have thus proven (\ref{sol.block2x2.tridet.indstep.minor.pf.8}) in each of
the two Subcases 2.1 and 2.2. Since these two Subcases cover the whole Case 2,
this shows that (\ref{sol.block2x2.tridet.indstep.minor.pf.8}) holds in Case
2.
\par
Thus, we have shown that (\ref{sol.block2x2.tridet.indstep.minor.pf.8}) holds
in each of the two Cases 1 and 2. Hence,
(\ref{sol.block2x2.tridet.indstep.minor.pf.8}) always holds. This completes
the proof of (\ref{sol.block2x2.tridet.indstep.minor.pf.8}).]
\par
Now, (\ref{sol.block2x2.tridet.indstep.minor.pf.2}) becomes%
\begin{align*}
U_{\sim\left(  n+M\right)  ,\sim\left(  n+q\right)  }  &  =\left(
\underbrace{u_{i,\mathbf{d}_{n+q}\left(  j\right)  }}_{\substack{=%
\begin{cases}
a_{i,j} & \text{if }i\leq n\text{ and }j\leq n;\\
b_{i,\mathbf{d}_{q}\left(  j-n\right)  }, & \text{if }i\leq n\text{ and
}j>n;\\
0, & \text{if }i>n\text{ and }j\leq n;\\
d_{i-n,\mathbf{d}_{q}\left(  j-n\right)  }, & \text{if }i>n\text{ and }j>n
\end{cases}
\\\text{(by (\ref{sol.block2x2.tridet.indstep.minor.pf.8}))}}}\right)  _{1\leq
i\leq n+M-1,\ 1\leq j\leq n+M-1}\\
&  =\left(
\begin{cases}
a_{i,j} & \text{if }i\leq n\text{ and }j\leq n;\\
b_{i,\mathbf{d}_{q}\left(  j-n\right)  }, & \text{if }i\leq n\text{ and
}j>n;\\
0, & \text{if }i>n\text{ and }j\leq n;\\
d_{i-n,\mathbf{d}_{q}\left(  j-n\right)  }, & \text{if }i>n\text{ and }j>n
\end{cases}
\right)  _{1\leq i\leq n+M-1,\ 1\leq j\leq n+M-1}\\
&  =\left(
\begin{array}
[c]{cc}%
A & B_{q}^{\prime}\\
0_{\left(  M-1\right)  \times n} & D_{\sim M,\sim q}%
\end{array}
\right)
\end{align*}
(by (\ref{sol.block2x2.tridet.indstep.minor.pf.6})). This proves
(\ref{sol.block2x2.tridet.indstep.minor}).}. Hence, for every $q\in\left\{
1,2,\ldots,M\right\}  $, we have%
\begin{equation}
\det\left(  U_{\sim\left(  n+M\right)  ,\sim\left(  n+q\right)  }\right)
=\det A\cdot\det\left(  D_{\sim M,\sim q}\right)
\label{sol.block2x2.tridet.indstep.minor-det}%
\end{equation}
\footnote{\textit{Proof of (\ref{sol.block2x2.tridet.indstep.minor-det}):} Let
$q\in\left\{  1,2,\ldots,M\right\}  $. We have assumed that Exercise
\ref{exe.block2x2.tridet} holds for $m=M-1$. Hence, we can apply Exercise
\ref{exe.block2x2.tridet} to $M-1$, $B_{q}^{\prime}$ and $D_{\sim M,\sim q}$
instead of $m$, $B$ and $D$. As a result, we obtain%
\[
\det\left(
\begin{array}
[c]{cc}%
A & B_{q}^{\prime}\\
0_{\left(  M-1\right)  \times n} & D_{\sim M,\sim q}%
\end{array}
\right)  =\det A\cdot\det\left(  D_{\sim M,\sim q}\right)  .
\]
Since $\left(
\begin{array}
[c]{cc}%
A & B_{q}^{\prime}\\
0_{\left(  M-1\right)  \times n} & D_{\sim M,\sim q}%
\end{array}
\right)  =U_{\sim\left(  n+M\right)  ,\sim\left(  n+q\right)  }$, this
rewrites as $\det\left(  U_{\sim\left(  n+M\right)  ,\sim\left(  n+q\right)
}\right)  =\det A\cdot\det\left(  D_{\sim M,\sim q}\right)  $. This proves
(\ref{sol.block2x2.tridet.indstep.minor-det}).}.

Furthermore, for every $q\in\left\{  1,2,\ldots,n\right\}  $, we have%
\begin{equation}
u_{n+M,q}=0 \label{sol.block2x2.tridet.indstep.u=0}%
\end{equation}
\footnote{\textit{Proof of (\ref{sol.block2x2.tridet.indstep.u=0}):} Let
$q\in\left\{  1,2,\ldots,n\right\}  $. Thus, $q\leq n$. Also,
$n+\underbrace{M}_{>0}>n$. Now, (\ref{sol.block2x2.tridet.indstep.minor.pf.1})
(applied to $\left(  i,j\right)  =\left(  n+M,q\right)  $) yields%
\[
u_{n+M,q}=%
\begin{cases}
a_{n+M,q} & \text{if }n+M\leq n\text{ and }q\leq n;\\
b_{n+M,q-n}, & \text{if }n+M\leq n\text{ and }q>n;\\
0, & \text{if }n+M>n\text{ and }q\leq n;\\
d_{n+M-n,q-n}, & \text{if }n+M>n\text{ and }q>n
\end{cases}
=0
\]
(since $n+M>n$ and $q\leq n$). This proves
(\ref{sol.block2x2.tridet.indstep.u=0}).}. On the other hand, for every
$q\in\left\{  1,2,\ldots,M\right\}  $, we have%
\begin{equation}
u_{n+M,n+q}=d_{M,q} \label{sol.block2x2.tridet.indstep.u=d}%
\end{equation}
\footnote{\textit{Proof of (\ref{sol.block2x2.tridet.indstep.u=d}):} Let
$q\in\left\{  1,2,\ldots,M\right\}  $. Thus, $q>0$, so that $n+\underbrace{q}%
_{>0}>n$. Also, $n+\underbrace{M}_{>0}>n$. Now,
(\ref{sol.block2x2.tridet.indstep.minor.pf.1}) (applied to $\left(
i,j\right)  =\left(  n+M,n+q\right)  $) yields%
\begin{align*}
u_{n+M,n+q}  &  =%
\begin{cases}
a_{n+M,n+q} & \text{if }n+M\leq n\text{ and }n+q\leq n;\\
b_{n+M,n+q-n}, & \text{if }n+M\leq n\text{ and }n+q>n;\\
0, & \text{if }n+M>n\text{ and }n+q\leq n;\\
d_{n+M-n,n+q-n}, & \text{if }n+M>n\text{ and }n+q>n
\end{cases}
\\
&  =d_{n+M-n,n+q-n}\ \ \ \ \ \ \ \ \ \ \left(  \text{since }n+M>n\text{ and
}n+q>n\right) \\
&  =d_{M,q}\ \ \ \ \ \ \ \ \ \ \left(  \text{since }n+M-n=M\text{ and
}n+q-n=q\right)  .
\end{align*}
This proves (\ref{sol.block2x2.tridet.indstep.u=d}).}.

Now, recall that $U=\left(  u_{i,j}\right)  _{1\leq i\leq n+M,\ 1\leq j\leq
n+M}$ and $n+M\in\left\{  1,2,\ldots,n+M\right\}  $ (since $n+M$ is positive).
Hence, Theorem \ref{thm.laplace.gen} \textbf{(a)} (applied to $n+M$, $U$,
$u_{i,j}$ and $n+M$ instead of $n$, $A$, $a_{i,j}$ and $p$) shows that%
\begin{align}
\det U  &  =\sum_{q=1}^{n+M}\left(  -1\right)  ^{\left(  n+M\right)
+q}u_{n+M,q}\det\left(  U_{\sim\left(  n+M\right)  ,\sim q}\right) \nonumber\\
&  =\sum_{q=1}^{n}\left(  -1\right)  ^{\left(  n+M\right)  +q}%
\underbrace{u_{n+M,q}}_{\substack{=0\\\text{(by
(\ref{sol.block2x2.tridet.indstep.u=0}))}}}\det\left(  U_{\sim\left(
n+M\right)  ,\sim q}\right) \nonumber\\
&  \ \ \ \ \ \ \ \ \ \ +\sum_{q=n+1}^{n+M}\left(  -1\right)  ^{\left(
n+M\right)  +q}u_{n+M,q}\det\left(  U_{\sim\left(  n+M\right)  ,\sim q}\right)
\nonumber\\
&  \ \ \ \ \ \ \ \ \ \ \left(  \text{since }0\leq n\leq n+M\right) \nonumber\\
&  =\underbrace{\sum_{q=1}^{n}\left(  -1\right)  ^{\left(  n+M\right)
+q}0\det\left(  U_{\sim\left(  n+M\right)  ,\sim q}\right)  }_{=0}\nonumber\\
&  \ \ \ \ \ \ \ \ \ \ +\sum_{q=n+1}^{n+M}\left(  -1\right)  ^{\left(
n+M\right)  +q}u_{n+M,q}\det\left(  U_{\sim\left(  n+M\right)  ,\sim q}\right)
\nonumber\\
&  =\sum_{q=n+1}^{n+M}\left(  -1\right)  ^{\left(  n+M\right)  +q}%
u_{n+M,q}\det\left(  U_{\sim\left(  n+M\right)  ,\sim q}\right) \nonumber\\
&  =\sum_{q=1}^{M}\underbrace{\left(  -1\right)  ^{\left(  n+M\right)
+\left(  n+q\right)  }}_{\substack{=\left(  -1\right)  ^{M+q}\\\text{(since
}\left(  n+M\right)  +\left(  n+q\right)  =2n+M+q\equiv M+q\operatorname{mod}%
2\text{)}}}\underbrace{u_{n+M,n+q}}_{\substack{=d_{M,q}\\\text{(by
(\ref{sol.block2x2.tridet.indstep.u=d}))}}}\underbrace{\det\left(
U_{\sim\left(  n+M\right)  ,\sim\left(  n+q\right)  }\right)  }%
_{\substack{=\det A\cdot\det\left(  D_{\sim M,\sim q}\right)  \\\text{(by
(\ref{sol.block2x2.tridet.indstep.minor-det}))}}}\nonumber\\
&  \ \ \ \ \ \ \ \ \ \ \left(  \text{here, we have substituted }n+q\text{ for
}q\text{ in the sum}\right) \nonumber\\
&  =\sum_{q=1}^{M}\left(  -1\right)  ^{M+q}d_{M,q}\det A\cdot\det\left(
D_{\sim M,\sim q}\right) \nonumber\\
&  =\det A\cdot\sum_{q=1}^{M}\left(  -1\right)  ^{M+q}d_{M,q}\det\left(
D_{\sim M,\sim q}\right)  . \label{sol.block2x2.tridet.indstep.detU}%
\end{align}

On the other hand, $D=\left(  d_{i,j}\right)  _{1\leq i\leq M,\ 1\leq j\leq
M}$ and $M\in\left\{  1,2,\ldots,M\right\}  $ (since $M>0$). Hence, Theorem
\ref{thm.laplace.gen} \textbf{(a)} (applied to $M$, $D$, $d_{i,j}$ and $M$
instead of $n$, $A$, $a_{i,j}$ and $p$) shows that%
\[
\det D=\sum_{q=1}^{M}\left(  -1\right)  ^{M+q}d_{M,q}\det\left(  D_{\sim
M,\sim q}\right)  .
\]
Thus, (\ref{sol.block2x2.tridet.indstep.detU}) becomes%
\[
\det U=\det A\cdot\underbrace{\sum_{q=1}^{M}\left(  -1\right)  ^{M+q}%
d_{M,q}\det\left(  D_{\sim M,\sim q}\right)  }_{=\det D}=\det A\cdot\det D.
\]
Since $U=\left(
\begin{array}
[c]{cc}%
A & B\\
0_{M\times n} & D
\end{array}
\right)  $, this rewrites as $\det\left(
\begin{array}
[c]{cc}%
A & B\\
0_{M\times n} & D
\end{array}
\right)  =\det A\cdot\det D$.

Now, let us forget that we fixed $n$, $A$, $B$ and $D$. We thus have shown
that for every $n\in\mathbb{N}$, for every $n\times n$-matrix $A$, for every
$n\times M$-matrix $B$, and for every $M\times M$-matrix $D$, we have
$\det\left(
\begin{array}
[c]{cc}%
A & B\\
0_{M\times n} & D
\end{array}
\right)  =\det A\cdot\det D$. In other words, Exercise
\ref{exe.block2x2.tridet} holds for $m=M$. This completes the induction step.
Hence, Exercise \ref{exe.block2x2.tridet} is solved by induction.
\end{proof}
\end{verlong}

\subsection{Solution to Exercise \ref{exe.block2x2.tridet.transposed}}

There are two ways to solve Exercise \ref{exe.block2x2.tridet.transposed}: One
way is to essentially repeat our above solution to Exercise
\ref{exe.block2x2.tridet} with some straightforward modifications (for
example, we must use Theorem \ref{thm.laplace.gen} \textbf{(b)} instead of
Theorem \ref{thm.laplace.gen} \textbf{(a)}). Another way is to derive Exercise
\ref{exe.block2x2.tridet.transposed} from Exercise \ref{exe.block2x2.tridet}
using transpose matrices. Let me show the latter way. We begin with a simple lemma:

\begin{lemma}
\label{lem.block2x2.transpose}Let $n$, $n^{\prime}$, $m$ and $m^{\prime}$ be
four nonnegative integers. Let $A\in\mathbb{K}^{n\times m}$, $B\in
\mathbb{K}^{n\times m^{\prime}}$, $C\in\mathbb{K}^{n^{\prime}\times m}$ and
$D\in\mathbb{K}^{n^{\prime}\times m^{\prime}}$. Then,%
\[
\left(
\begin{array}
[c]{cc}%
A & B\\
C & D
\end{array}
\right)  ^{T}=\left(
\begin{array}
[c]{cc}%
A^{T} & C^{T}\\
B^{T} & D^{T}%
\end{array}
\right)  .
\]

\end{lemma}

\begin{vershort}
\begin{proof}
[Proof of Lemma \ref{lem.block2x2.transpose}.]Lemma
\ref{lem.block2x2.transpose} results in a straightforward way by recalling the
definitions of $\left(
\begin{array}
[c]{cc}%
A & B\\
C & D
\end{array}
\right)  ^{T}$ and $\left(
\begin{array}
[c]{cc}%
A^{T} & C^{T}\\
B^{T} & D^{T}%
\end{array}
\right)  $ and comparing.
\end{proof}
\end{vershort}

\begin{verlong}
\begin{proof}
[Proof of Lemma \ref{lem.block2x2.transpose}.]Write the $n\times m$-matrix $A$
in the form $A=\left(  a_{i,j}\right)  _{1\leq i\leq n,\ 1\leq j\leq m}$.
Thus, $A^{T}=\left(  a_{j,i}\right)  _{1\leq i\leq m,\ 1\leq j\leq n}$ (by the
definition of $A^{T}$).

Write the $n\times m^{\prime}$-matrix $B$ in the form $B=\left(
b_{i,j}\right)  _{1\leq i\leq n,\ 1\leq j\leq m^{\prime}}$. Thus,
$B^{T}=\left(  b_{j,i}\right)  _{1\leq i\leq m^{\prime},\ 1\leq j\leq n}$ (by
the definition of $B^{T}$).

Write the $n^{\prime}\times m$-matrix $C$ in the form $C=\left(
c_{i,j}\right)  _{1\leq i\leq n^{\prime},\ 1\leq j\leq m}$. Thus,
$C^{T}=\left(  c_{j,i}\right)  _{1\leq i\leq m,\ 1\leq j\leq n^{\prime}}$ (by
the definition of $C^{T}$).

Write the $n^{\prime}\times m^{\prime}$-matrix $D$ in the form $D=\left(
d_{i,j}\right)  _{1\leq i\leq n^{\prime},\ 1\leq j\leq m^{\prime}}$. Thus,
$D^{T}=\left(  d_{j,i}\right)  _{1\leq i\leq m^{\prime},\ 1\leq j\leq
n^{\prime}}$ (by the definition of $D^{T}$).

We have $A=\left(  a_{i,j}\right)  _{1\leq i\leq n,\ 1\leq j\leq m}$,
$B=\left(  b_{i,j}\right)  _{1\leq i\leq n,\ 1\leq j\leq m^{\prime}}$,
$C=\left(  c_{i,j}\right)  _{1\leq i\leq n^{\prime},\ 1\leq j\leq m}$ and
$D=\left(  d_{i,j}\right)  _{1\leq i\leq n^{\prime},\ 1\leq j\leq m^{\prime}}%
$. Hence, (\ref{eq.def.block2x2.formal}) yields%
\[
\left(
\begin{array}
[c]{cc}%
A & B\\
C & D
\end{array}
\right)  =\left(
\begin{cases}
a_{i,j}, & \text{if }i\leq n\text{ and }j\leq m;\\
b_{i,j-m}, & \text{if }i\leq n\text{ and }j>m;\\
c_{i-n,j}, & \text{if }i>n\text{ and }j\leq m;\\
d_{i-n,j-m}, & \text{if }i>n\text{ and }j>m
\end{cases}
\right)  _{1\leq i\leq n+n^{\prime},\ 1\leq j\leq m+m^{\prime}}.
\]
Hence, the definition of the transpose matrix $\left(
\begin{array}
[c]{cc}%
A & B\\
C & D
\end{array}
\right)  ^{T}$ yields%
\begin{align}
&  \left(
\begin{array}
[c]{cc}%
A & B\\
C & D
\end{array}
\right)  ^{T}\nonumber\\
&  =\left(
\begin{cases}
a_{j,i}, & \text{if }j\leq n\text{ and }i\leq m;\\
b_{j,i-m}, & \text{if }j\leq n\text{ and }i>m;\\
c_{j-n,i}, & \text{if }j>n\text{ and }i\leq m;\\
d_{j-n,i-m}, & \text{if }j>n\text{ and }i>m
\end{cases}
\right)  _{1\leq i\leq m+m^{\prime},\ 1\leq j\leq n+n^{\prime}}.
\label{pf.lem.block2x2.transpose.1}%
\end{align}

On the other hand, $A^{T}=\left(  a_{j,i}\right)  _{1\leq i\leq m,\ 1\leq
j\leq n}$, $C^{T}=\left(  c_{j,i}\right)  _{1\leq i\leq m,\ 1\leq j\leq
n^{\prime}}$, $B^{T}=\left(  b_{j,i}\right)  _{1\leq i\leq m^{\prime},\ 1\leq
j\leq n}$ and $D^{T}=\left(  d_{j,i}\right)  _{1\leq i\leq m^{\prime},\ 1\leq
j\leq n^{\prime}}$. Hence, (\ref{eq.def.block2x2.formal}) (applied to $m$,
$m^{\prime}$, $n$, $n^{\prime}$, $A^{T}$, $C^{T}$, $B^{T}$, $D^{T}$, $a_{j,i}%
$, $c_{j,i}$, $b_{j,i}$ and $d_{j,i}$ instead of $n$, $n^{\prime}$, $m$,
$m^{\prime}$, $A$, $B$, $C$, $D$, $a_{i,j}$, $b_{i,j}$, $c_{i,j}$ and
$d_{i,j}$) yields%
\begin{align}
&  \left(
\begin{array}
[c]{cc}%
A^{T} & C^{T}\\
B^{T} & D^{T}%
\end{array}
\right) \nonumber\\
&  =\left(
\begin{cases}
a_{j,i}, & \text{if }i\leq m\text{ and }j\leq n;\\
c_{j-n,i}, & \text{if }i\leq m\text{ and }j>n;\\
b_{j,i-m}, & \text{if }i>m\text{ and }j\leq n;\\
d_{j-n,i-m}, & \text{if }i>m\text{ and }j>n
\end{cases}
\right)  _{1\leq i\leq m+m^{\prime},\ 1\leq j\leq n+n^{\prime}}.
\label{pf.lem.block2x2.transpose.2}%
\end{align}

But every $\left(  i,j\right)  \in\left\{  1,2,\ldots,m+m^{\prime}\right\}
\times\left\{  1,2,\ldots,n+n^{\prime}\right\}  $ satisfies%
\begin{align*}
&
\begin{cases}
a_{j,i}, & \text{if }j\leq n\text{ and }i\leq m;\\
b_{j,i-m}, & \text{if }j\leq n\text{ and }i>m;\\
c_{j-n,i}, & \text{if }j>n\text{ and }i\leq m;\\
d_{j-n,i-m}, & \text{if }j>n\text{ and }i>m
\end{cases}
\\
&  =%
\begin{cases}
a_{j,i}, & \text{if }j\leq n\text{ and }i\leq m;\\
c_{j-n,i}, & \text{if }j>n\text{ and }i\leq m;\\
b_{j,i-m}, & \text{if }j\leq n\text{ and }i>m;\\
d_{j-n,i-m}, & \text{if }j>n\text{ and }i>m
\end{cases}
\\
&  \ \ \ \ \ \ \ \ \ \ \left(  \text{here, we have swapped the second and the
third cases}\right) \\
&  =%
\begin{cases}
a_{j,i}, & \text{if }i\leq m\text{ and }j\leq n;\\
c_{j-n,i}, & \text{if }j>n\text{ and }i\leq m;\\
b_{j,i-m}, & \text{if }j\leq n\text{ and }i>m;\\
d_{j-n,i-m}, & \text{if }j>n\text{ and }i>m
\end{cases}
\\
&  \ \ \ \ \ \ \ \ \ \ \left(  \text{since }\left(  j\leq n\text{ and }i\leq
m\right)  \text{ is equivalent to }\left(  i\leq m\text{ and }j\leq n\right)
\right) \\
&  =%
\begin{cases}
a_{j,i}, & \text{if }i\leq m\text{ and }j\leq n;\\
c_{j-n,i}, & \text{if }i\leq m\text{ and }j>n;\\
b_{j,i-m}, & \text{if }j\leq n\text{ and }i>m;\\
d_{j-n,i-m}, & \text{if }j>n\text{ and }i>m
\end{cases}
\\
&  \ \ \ \ \ \ \ \ \ \ \left(  \text{since }\left(  j>n\text{ and }i\leq
m\right)  \text{ is equivalent to }\left(  i\leq m\text{ and }j>n\right)
\right)
\end{align*}%
\begin{align*}
&  =%
\begin{cases}
a_{j,i}, & \text{if }i\leq m\text{ and }j\leq n;\\
c_{j-n,i}, & \text{if }i\leq m\text{ and }j>n;\\
b_{j,i-m}, & \text{if }i>m\text{ and }j\leq n;\\
d_{j-n,i-m}, & \text{if }j>n\text{ and }i>m
\end{cases}
\\
&  \ \ \ \ \ \ \ \ \ \ \left(  \text{since }\left(  j\leq n\text{ and
}i>m\right)  \text{ is equivalent to }\left(  i>m\text{ and }j\leq n\right)
\right) \\
&  =%
\begin{cases}
a_{j,i}, & \text{if }i\leq m\text{ and }j\leq n;\\
c_{j-n,i}, & \text{if }i\leq m\text{ and }j>n;\\
b_{j,i-m}, & \text{if }i>m\text{ and }j\leq n;\\
d_{j-n,i-m}, & \text{if }i>m\text{ and }j>n
\end{cases}
\\
&  \ \ \ \ \ \ \ \ \ \ \left(  \text{since }\left(  j>n\text{ and }i>m\right)
\text{ is equivalent to }\left(  i>m\text{ and }j>n\right)  \right)  .
\end{align*}
In other words,%
\begin{align*}
&  \left(
\begin{cases}
a_{j,i}, & \text{if }j\leq n\text{ and }i\leq m;\\
b_{j,i-m}, & \text{if }j\leq n\text{ and }i>m;\\
c_{j-n,i}, & \text{if }j>n\text{ and }i\leq m;\\
d_{j-n,i-m}, & \text{if }j>n\text{ and }i>m
\end{cases}
\right)  _{1\leq i\leq m+m^{\prime},\ 1\leq j\leq n+n^{\prime}}\\
&  =\left(
\begin{cases}
a_{j,i}, & \text{if }i\leq m\text{ and }j\leq n;\\
c_{j-n,i}, & \text{if }i\leq m\text{ and }j>n;\\
b_{j,i-m}, & \text{if }i>m\text{ and }j\leq n;\\
d_{j-n,i-m}, & \text{if }i>m\text{ and }j>n
\end{cases}
\right)  _{1\leq i\leq m+m^{\prime},\ 1\leq j\leq n+n^{\prime}}.
\end{align*}
Comparing this with (\ref{pf.lem.block2x2.transpose.2}), we obtain%
\begin{align*}
\left(
\begin{array}
[c]{cc}%
A^{T} & C^{T}\\
B^{T} & D^{T}%
\end{array}
\right)   &  =\left(
\begin{cases}
a_{j,i}, & \text{if }j\leq n\text{ and }i\leq m;\\
b_{j,i-m}, & \text{if }j\leq n\text{ and }i>m;\\
c_{j-n,i}, & \text{if }j>n\text{ and }i\leq m;\\
d_{j-n,i-m}, & \text{if }j>n\text{ and }i>m
\end{cases}
\right)  _{1\leq i\leq m+m^{\prime},\ 1\leq j\leq n+n^{\prime}}\\
&  =\left(
\begin{array}
[c]{cc}%
A & B\\
C & D
\end{array}
\right)  ^{T}%
\end{align*}
(by (\ref{pf.lem.block2x2.transpose.1})). This proves Lemma
\ref{lem.block2x2.transpose}.
\end{proof}
\end{verlong}

Now, we can comfortably solve Exercise \ref{exe.block2x2.tridet.transposed}:

\begin{vershort}
\begin{proof}
[Solution to Exercise \ref{exe.block2x2.tridet.transposed}.]Exercise
\ref{exe.ps4.4} (applied to $n+m$ and $\left(
\begin{array}
[c]{cc}%
A & 0_{n\times m}\\
C & D
\end{array}
\right)  $ instead of $n$ and $A$) shows that
\[
\det\left(  \left(
\begin{array}
[c]{cc}%
A & 0_{n\times m}\\
C & D
\end{array}
\right)  ^{T}\right)  =\det\left(
\begin{array}
[c]{cc}%
A & 0_{n\times m}\\
C & D
\end{array}
\right)  .
\]
Hence,%
\begin{align*}
&  \det\left(
\begin{array}
[c]{cc}%
A & 0_{n\times m}\\
C & D
\end{array}
\right) \\
&  =\det\left(  \underbrace{\left(
\begin{array}
[c]{cc}%
A & 0_{n\times m}\\
C & D
\end{array}
\right)  ^{T}}_{\substack{=\left(
\begin{array}
[c]{cc}%
A^{T} & C^{T}\\
\left(  0_{n\times m}\right)  ^{T} & D^{T}%
\end{array}
\right)  \\\text{(by Lemma \ref{lem.block2x2.transpose} (applied to }n\text{,
}m\text{, }n\text{, }m\\\text{and }0_{n\times m}\text{ instead of }n\text{,
}n^{\prime}\text{, }m\text{, }m^{\prime}\text{ and }B\text{))}}}\right)
=\det\underbrace{\left(
\begin{array}
[c]{cc}%
A^{T} & C^{T}\\
\left(  0_{n\times m}\right)  ^{T} & D^{T}%
\end{array}
\right)  }_{\substack{=\left(
\begin{array}
[c]{cc}%
A^{T} & C^{T}\\
0_{m\times n} & D^{T}%
\end{array}
\right)  \\\text{(since }\left(  0_{n\times m}\right)  ^{T}=0_{m\times
n}\text{)}}}\\
&  =\det\left(
\begin{array}
[c]{cc}%
A^{T} & C^{T}\\
0_{m\times n} & D^{T}%
\end{array}
\right)  =\underbrace{\det\left(  A^{T}\right)  }_{\substack{=\det
A\\\text{(by Exercise \ref{exe.ps4.4})}}}\cdot\underbrace{\det\left(
D^{T}\right)  }_{\substack{=\det D\\\text{(by Exercise \ref{exe.ps4.4}%
}\\\text{(applied to }m\text{ and }D\text{ instead}\\\text{of }n\text{ and
}A\text{))}}}\\
&  \ \ \ \ \ \ \ \ \ \ \left(  \text{by Exercise \ref{exe.block2x2.tridet}
(applied to }A^{T}\text{, }C^{T}\text{ and }D^{T}\text{ instead of }A\text{,
}B\text{ and }D\text{)}\right) \\
&  =\det A\cdot\det D.
\end{align*}
This solves Exercise \ref{exe.block2x2.tridet.transposed}.
\end{proof}
\end{vershort}

\begin{verlong}
\begin{proof}
[Solution to Exercise \ref{exe.block2x2.tridet.transposed}.]The matrix $A$ is
an $n\times n$-matrix; hence, its transpose $A^{T}$ is an $n\times n$-matrix.
In other words, $A^{T}\in\mathbb{K}^{n\times n}$.

The matrix $C$ is an $m\times n$-matrix; hence, its transpose $C^{T}$ is an
$n\times m$-matrix. In other words, $C^{T}\in\mathbb{K}^{n\times m}$.

The matrix $D$ is an $m\times m$-matrix; hence, its transpose $D^{T}$ is an
$m\times m$-matrix. In other words, $D^{T}\in\mathbb{K}^{m\times m}$.

Clearly, $\left(  0_{n\times m}\right)  ^{T}=0_{m\times n}$.

The matrix $\left(
\begin{array}
[c]{cc}%
A & 0_{n\times m}\\
C & D
\end{array}
\right)  $ is an $\left(  n+m\right)  \times\left(  n+m\right)  $-matrix.
Hence, Exercise \ref{exe.ps4.4} (applied to $n+m$ and $\left(
\begin{array}
[c]{cc}%
A & 0_{n\times m}\\
C & D
\end{array}
\right)  $ instead of $n$ and $A$) shows that
\[
\det\left(  \left(
\begin{array}
[c]{cc}%
A & 0_{n\times m}\\
C & D
\end{array}
\right)  ^{T}\right)  =\det\left(
\begin{array}
[c]{cc}%
A & 0_{n\times m}\\
C & D
\end{array}
\right)  .
\]
Hence,%
\begin{align*}
&  \det\left(
\begin{array}
[c]{cc}%
A & 0_{n\times m}\\
C & D
\end{array}
\right) \\
&  =\det\left(  \underbrace{\left(
\begin{array}
[c]{cc}%
A & 0_{n\times m}\\
C & D
\end{array}
\right)  ^{T}}_{\substack{=\left(
\begin{array}
[c]{cc}%
A^{T} & C^{T}\\
\left(  0_{n\times m}\right)  ^{T} & D^{T}%
\end{array}
\right)  \\\text{(by Lemma \ref{lem.block2x2.transpose} (applied to }n\text{,
}m\text{, }n\text{, }m\\\text{and }0_{n\times m}\text{ instead of }n\text{,
}n^{\prime}\text{, }m\text{, }m^{\prime}\text{ and }B\text{))}}}\right)
=\det\underbrace{\left(
\begin{array}
[c]{cc}%
A^{T} & C^{T}\\
\left(  0_{n\times m}\right)  ^{T} & D^{T}%
\end{array}
\right)  }_{\substack{=\left(
\begin{array}
[c]{cc}%
A^{T} & C^{T}\\
0_{m\times n} & D^{T}%
\end{array}
\right)  \\\text{(since }\left(  0_{n\times m}\right)  ^{T}=0_{m\times
n}\text{)}}}\\
&  =\det\left(
\begin{array}
[c]{cc}%
A^{T} & C^{T}\\
0_{m\times n} & D^{T}%
\end{array}
\right)  =\underbrace{\det\left(  A^{T}\right)  }_{\substack{=\det
A\\\text{(by Exercise \ref{exe.ps4.4})}}}\cdot\underbrace{\det\left(
D^{T}\right)  }_{\substack{=\det D\\\text{(by Exercise \ref{exe.ps4.4}%
}\\\text{(applied to }m\text{ and }D\text{ instead}\\\text{of }n\text{ and
}A\text{))}}}\\
&  \ \ \ \ \ \ \ \ \ \ \left(  \text{by Exercise \ref{exe.block2x2.tridet}
(applied to }A^{T}\text{, }C^{T}\text{ and }D^{T}\text{ instead of }A\text{,
}B\text{ and }D\text{)}\right) \\
&  =\det A\cdot\det D.
\end{align*}
This solves Exercise \ref{exe.block2x2.tridet.transposed}.
\end{proof}
\end{verlong}

\subsection{Second solution to Exercise \ref{exe.ps4.5}}

\begin{proof}
[Second solution to Exercise \ref{exe.ps4.5} (sketched).]\textbf{(b)} We have%
\begin{align*}
&  \det\left(
\begin{array}
[c]{ccccc}%
a & b & c & d & e\\
f & 0 & 0 & 0 & g\\
h & 0 & 0 & 0 & i\\
j & 0 & 0 & 0 & k\\
\ell & m & n & o & p
\end{array}
\right) \\
&  =-\det\left(
\begin{array}
[c]{ccccc}%
a & b & c & d & e\\
\ell & m & n & o & p\\
h & 0 & 0 & 0 & i\\
j & 0 & 0 & 0 & k\\
f & 0 & 0 & 0 & g
\end{array}
\right) \\
&  \ \ \ \ \ \ \ \ \ \ \left(
\begin{array}
[c]{c}%
\text{by Exercise \ref{exe.ps4.6} \textbf{(a)}, because we have just}\\
\text{swapped the }2\text{-nd and the }5\text{-th rows of the matrix}%
\end{array}
\right) \\
&  =\det\left(
\begin{array}
[c]{ccccc}%
d & b & c & a & e\\
o & m & n & \ell & p\\
0 & 0 & 0 & h & i\\
0 & 0 & 0 & j & k\\
0 & 0 & 0 & f & g
\end{array}
\right) \\
&  \ \ \ \ \ \ \ \ \ \ \left(
\begin{array}
[c]{c}%
\text{by Exercise \ref{exe.ps4.6} \textbf{(b)}, because we have just}\\
\text{swapped the }1\text{-st and the }4\text{-th columns of the matrix}%
\end{array}
\right) \\
&  =\underbrace{\det\left(
\begin{array}
[c]{ccc}%
d & b & c\\
o & m & n\\
0 & 0 & 0
\end{array}
\right)  }_{\substack{=0\\\text{(by Exercise \ref{exe.ps4.6} \textbf{(c)})}%
}}\cdot\det\left(
\begin{array}
[c]{cc}%
j & k\\
f & g
\end{array}
\right) \\
&  \ \ \ \ \ \ \ \ \ \ \left(
\begin{array}
[c]{c}%
\text{by Exercise \ref{exe.block2x2.tridet}, applied to }\left(
\begin{array}
[c]{ccc}%
d & b & c\\
o & m & n\\
0 & 0 & 0
\end{array}
\right)  \text{,}\\
\left(
\begin{array}
[c]{cc}%
a & e\\
\ell & p\\
h & i
\end{array}
\right)  \text{, }\left(
\begin{array}
[c]{cc}%
j & k\\
f & g
\end{array}
\right)  \text{, }2\text{ and }3\text{ instead of }A\text{, }B\text{,
}D\text{, }m\text{ and }n
\end{array}
\right) \\
&  =0.
\end{align*}
This solves Exercise \ref{exe.ps4.5} \textbf{(b)}.

\textbf{(a)} We have%
\begin{align*}
&  \det\left(
\begin{array}
[c]{cccc}%
a & b & c & d\\
\ell & 0 & 0 & e\\
k & 0 & 0 & f\\
j & i & h & g
\end{array}
\right) \\
&  =-\det\left(
\begin{array}
[c]{cccc}%
a & b & c & d\\
j & i & h & g\\
k & 0 & 0 & f\\
\ell & 0 & 0 & e
\end{array}
\right) \\
&  \ \ \ \ \ \ \ \ \ \ \left(
\begin{array}
[c]{c}%
\text{by Exercise \ref{exe.ps4.6} \textbf{(a)}, because we have just}\\
\text{swapped the }2\text{-nd and the }4\text{-th rows of the matrix}%
\end{array}
\right) \\
&  =\det\left(
\begin{array}
[c]{cccc}%
c & b & a & d\\
h & i & j & g\\
0 & 0 & k & f\\
0 & 0 & \ell & e
\end{array}
\right) \\
&  \ \ \ \ \ \ \ \ \ \ \left(
\begin{array}
[c]{c}%
\text{by Exercise \ref{exe.ps4.6} \textbf{(b)}, because we have just}\\
\text{swapped the }1\text{-st and the }3\text{-th columns of the matrix}%
\end{array}
\right) \\
&  =\underbrace{\det\left(
\begin{array}
[c]{cc}%
c & b\\
h & i
\end{array}
\right)  }_{=ci-bh}\cdot\underbrace{\det\left(
\begin{array}
[c]{cc}%
k & f\\
\ell & e
\end{array}
\right)  }_{=ek-\ell f}\\
&  \ \ \ \ \ \ \ \ \ \ \left(
\begin{array}
[c]{c}%
\text{by Exercise \ref{exe.block2x2.tridet}, applied to }\left(
\begin{array}
[c]{cc}%
c & b\\
h & i
\end{array}
\right)  \text{,}\\
\left(
\begin{array}
[c]{cc}%
a & d\\
j & g
\end{array}
\right)  \text{, }\left(
\begin{array}
[c]{cc}%
k & f\\
\ell & e
\end{array}
\right)  \text{, }2\text{ and }2\text{ instead of }A\text{, }B\text{,
}D\text{, }m\text{ and }n
\end{array}
\right) \\
&  =\left(  ci-bh\right)  \left(  ek-\ell f\right)  =\left(  bh-ci\right)
\left(  \ell f-ek\right)  .
\end{align*}
This solves Exercise \ref{exe.ps4.5} \textbf{(a)}. (Notice that we have
obtained the result in its factored form!)
\end{proof}

\subsection{\label{sect.sol.det.anotherpattern}Solution to Exercise
\ref{exe.det.anotherpattern}}

Before we solve Exercise \ref{exe.det.anotherpattern}, let us introduce some
notation that enables us to clearly speak about permutations of rows and
columns in a matrix:

\begin{definition}
\label{def.sol.addexe.jacobi-complement.Agd}Let $n\in\mathbb{N}$ and
$m\in\mathbb{N}$. Let $A=\left(  a_{i,j}\right)  _{1\leq i\leq n,\ 1\leq j\leq
m}\in\mathbb{K}^{n\times m}$. Let $\gamma\in S_{n}$ and $\delta\in S_{m}$.
Then, $A_{\left[  \gamma,\delta\right]  }$ denotes the matrix $\left(
a_{\gamma\left(  i\right)  ,\delta\left(  j\right)  }\right)  _{1\leq i\leq
n,\ 1\leq j\leq m}\in\mathbb{K}^{n\times m}$.
\end{definition}

We shall use this definition throughout Section
\ref{sect.sol.det.anotherpattern}.

\begin{example}
For this example, let $n=3$ and $m=4$, and let $A=\left(
\begin{array}
[c]{cccc}%
a & b & c & d\\
a^{\prime} & b^{\prime} & c^{\prime} & d^{\prime}\\
a^{\prime\prime} & b^{\prime\prime} & c^{\prime\prime} & d^{\prime\prime}%
\end{array}
\right)  \in\mathbb{K}^{3\times4}$ be a $3\times4$-matrix. Let $\gamma\in
S_{3}$ be the permutation that sends $1,2,3$ to $2,3,1$, respectively. Let
$\delta\in S_{4}$ be the permutation that sends $1,2,3,4$ to $2,1,4,3$,
respectively. Then,%
\[
A_{\left[  \gamma,\delta\right]  }=\left(
\begin{array}
[c]{cccc}%
b^{\prime} & a^{\prime} & d^{\prime} & c^{\prime}\\
b^{\prime\prime} & a^{\prime\prime} & d^{\prime\prime} & c^{\prime\prime}\\
b & a & d & c
\end{array}
\right)  .
\]

\end{example}

\begin{remark}
Let $n$, $m$, $A$, $\gamma$ and $\delta$ be as in Definition
\ref{def.sol.addexe.jacobi-complement.Agd}. Then, it is easy to see that
$A_{\left[  \gamma,\delta\right]  }=\operatorname*{sub}\nolimits_{\gamma
\left(  1\right)  ,\gamma\left(  2\right)  ,\ldots,\gamma\left(  n\right)
}^{\delta\left(  1\right)  ,\delta\left(  2\right)  ,\ldots,\delta\left(
m\right)  }A$ (where we are using the notation introduced in Definition
\ref{def.submatrix}). Visually speaking, $A_{\left[  \gamma,\delta\right]  }$
is the matrix obtained from $A$ by permuting the rows (using the permutation
$\gamma$) and permuting the columns (using the permutation $\delta$).
\end{remark}

\begin{remark}
Let $n\in\mathbb{N}$, and let $B\in\mathbb{K}^{n\times n}$. Let $\kappa\in
S_{n}$. Then, Definition \ref{def.sol.addexe.jacobi-complement.Agd} gives rise
to an $n\times n$-matrix $B_{\left[  \kappa,\operatorname*{id}\right]  }$
(where $\operatorname*{id}$ denotes the identity permutation
$\operatorname*{id}\nolimits_{\left\{  1,2,\ldots,n\right\}  }\in S_{n}$).
This matrix $B_{\left[  \kappa,\operatorname*{id}\right]  }$ is precisely the
matrix $B_{\kappa}$ from Lemma \ref{lem.det.sigma}.
\end{remark}

The following lemma will be crucial for us:

\begin{lemma}
\label{lem.sol.addexe.jacobi-complement.Agd.det}Let $n\in\mathbb{N}$. Let
$A\in\mathbb{K}^{n\times n}$. Let $\gamma\in S_{n}$ and $\delta\in S_{n}$.
Then,%
\[
\det\left(  A_{\left[  \gamma,\delta\right]  }\right)  =\left(  -1\right)
^{\gamma}\left(  -1\right)  ^{\delta}\det A.
\]

\end{lemma}

Lemma \ref{lem.sol.addexe.jacobi-complement.Agd.det} simply says that if we
permute the rows and permute the columns of a square matrix, then the
determinant of this matrix gets multiplied by the product of the signs of the
two permutations.

\begin{proof}
[Proof of Lemma \ref{lem.sol.addexe.jacobi-complement.Agd.det}.]Let $\left[
n\right]  $ denote the set $\left\{  1,2,\ldots,n\right\}  $.

\begin{verlong}
We have $\gamma\in S_{n}$. In other words, $\gamma$ is a permutation of
$\left\{  1,2,\ldots,n\right\}  $ (since $S_{n}$ is the set of all
permutations of $\left\{  1,2,\ldots,n\right\}  $). In other words, $\gamma$
is a permutation of $\left[  n\right]  $ (since $\left[  n\right]  =\left\{
1,2,\ldots,n\right\}  $). In other words, $\gamma$ is a bijective map $\left[
n\right]  \rightarrow\left[  n\right]  $. The same argument (but applied to
$\delta$ instead of $\gamma$) shows that $\delta$ is a bijective map $\left[
n\right]  \rightarrow\left[  n\right]  $.
\end{verlong}

Write the $n\times n$-matrix $A$ in the form $A=\left(  a_{i,j}\right)
_{1\leq i\leq n,\ 1\leq j\leq n}$. Then, \newline$A^{T}=\left(  a_{j,i}%
\right)  _{1\leq i\leq n,\ 1\leq j\leq n}$ (by the definition of $A^{T}$).

Define a new matrix $C$ by $C=\left(  a_{j,\delta\left(  i\right)  }\right)
_{1\leq i\leq n,\ 1\leq j\leq n}$. Then, Lemma \ref{lem.det.sigma}
\textbf{(a)} (applied to $\delta$, $A^{T}$, $a_{j,i}$ and $C$ instead of
$\kappa$, $B$, $b_{i,j}$ and $B_{\kappa}$) yields $\det C=\left(  -1\right)
^{\delta}\cdot\det\left(  A^{T}\right)  $ (because $A^{T}=\left(
a_{j,i}\right)  _{1\leq i\leq n,\ 1\leq j\leq n}$ and $C=\left(
a_{j,\delta\left(  i\right)  }\right)  _{1\leq i\leq n,\ 1\leq j\leq n}$).
Thus,%
\[
\det C=\left(  -1\right)  ^{\delta}\cdot\underbrace{\det\left(  A^{T}\right)
}_{\substack{=\det A\\\text{(by Exercise \ref{exe.ps4.4})}}}=\left(
-1\right)  ^{\delta}\cdot\det A.
\]

Clearly, $C$ is an $n\times n$-matrix. Hence, Exercise \ref{exe.ps4.4}
(applied to $C$ instead of $A$) yields
\[
\det\left(  C^{T}\right)  =\det C=\left(  -1\right)  ^{\delta}\cdot\det A.
\]

But we have $C=\left(  a_{j,\delta\left(  i\right)  }\right)  _{1\leq i\leq
n,\ 1\leq j\leq n}$. Hence, $C^{T}=\left(  a_{i,\delta\left(  j\right)
}\right)  _{1\leq i\leq n,\ 1\leq j\leq n}$ (by the definition of $C^{T}$). We
have $C^{T}=\left(  a_{i,\delta\left(  j\right)  }\right)  _{1\leq i\leq
n,\ 1\leq j\leq n}$ and $A_{\left[  \gamma,\delta\right]  }=\left(
a_{\gamma\left(  i\right)  ,\delta\left(  j\right)  }\right)  _{1\leq i\leq
n,\ 1\leq j\leq n}$ (by the definition of $A_{\left[  \gamma,\delta\right]  }%
$). Thus, Lemma \ref{lem.det.sigma} \textbf{(a)} (applied to $\gamma$, $C^{T}%
$, $a_{i,\delta\left(  j\right)  }$ and $A_{\left[  \gamma,\delta\right]  }$
instead of $\kappa$, $B$, $b_{i,j}$ and $B_{\kappa}$) yields%
\[
\det\left(  A_{\left[  \gamma,\delta\right]  }\right)  =\left(  -1\right)
^{\gamma}\cdot\underbrace{\det\left(  C^{T}\right)  }_{=\left(  -1\right)
^{\delta}\cdot\det A}=\left(  -1\right)  ^{\gamma}\cdot\left(  -1\right)
^{\delta}\cdot\det A=\left(  -1\right)  ^{\gamma}\left(  -1\right)  ^{\delta
}\det A.
\]
This proves Lemma \ref{lem.sol.addexe.jacobi-complement.Agd.det}.
\end{proof}

\begin{corollary}
\label{cor.sol.det.anotherpattern.gammagamma}Let $n\in\mathbb{N}$. Let
$A\in\mathbb{K}^{n\times n}$. Let $\gamma\in S_{n}$. Then, $\det\left(
A_{\left[  \gamma,\gamma\right]  }\right)  =\det A$.
\end{corollary}

\begin{proof}
[Proof of Corollary \ref{cor.sol.det.anotherpattern.gammagamma}.]We have
$\left(  -1\right)  ^{\gamma}=\left(  -1\right)  ^{\ell\left(  \gamma\right)
}$ (by the definition of $\left(  -1\right)  ^{\gamma}$). Hence,%
\[
\left(  \underbrace{\left(  -1\right)  ^{\gamma}}_{=\left(  -1\right)
^{\ell\left(  \gamma\right)  }}\right)  ^{2}=\left(  \left(  -1\right)
^{\ell\left(  \gamma\right)  }\right)  ^{2}=\left(  -1\right)  ^{\ell\left(
\gamma\right)  \cdot2}=1
\]
(since the integer $\ell\left(  \gamma\right)  \cdot2$ is even). But Lemma
\ref{lem.sol.addexe.jacobi-complement.Agd.det} (applied to $\delta=\gamma$)
yields%
\[
\det\left(  A_{\left[  \gamma,\gamma\right]  }\right)  =\underbrace{\left(
-1\right)  ^{\gamma}\left(  -1\right)  ^{\gamma}}_{=\left(  \left(  -1\right)
^{\gamma}\right)  ^{2}=1}\det A=\det A.
\]
This proves Corollary \ref{cor.sol.det.anotherpattern.gammagamma}.
\end{proof}

\begin{proof}
[Solution to Exercise \ref{exe.det.anotherpattern}.]\textbf{(a)} Let
$a,a^{\prime},a^{\prime\prime},b,b^{\prime},b^{\prime\prime},c,c^{\prime
},c^{\prime\prime},d,d^{\prime},d^{\prime\prime},e$ be elements of
$\mathbb{K}$. Let $A$ be the $7\times7$-matrix%
\[
\left(
\begin{array}
[c]{ccccccc}%
a & 0 & 0 & 0 & 0 & 0 & b\\
0 & a^{\prime} & 0 & 0 & 0 & b^{\prime} & 0\\
0 & 0 & a^{\prime\prime} & 0 & b^{\prime\prime} & 0 & 0\\
0 & 0 & 0 & e & 0 & 0 & 0\\
0 & 0 & c^{\prime\prime} & 0 & d^{\prime\prime} & 0 & 0\\
0 & c^{\prime} & 0 & 0 & 0 & d^{\prime} & 0\\
c & 0 & 0 & 0 & 0 & 0 & d
\end{array}
\right)  .
\]
Thus, the exercise demands that we find $\det A$.

Let $\gamma\in S_{7}$ be the permutation of the set $\left\{
1,2,3,4,5,6,7\right\}  $ which sends $1,2,3,4,5,6,7$ to $1,7,3,5,4,6,2$,
respectively. (This permutation $\gamma$ swaps $2$ with $7$ and swaps $4$ with
$5$; it leaves $1,3,6$ unchanged.) Then, straightforward computation (using
the definitions of $A$ and of $A_{\left[  \gamma,\gamma\right]  }$) shows that%
\begin{equation}
A_{\left[  \gamma,\gamma\right]  }=\left(
\begin{array}
[c]{ccccccc}%
a & b & 0 & 0 & 0 & 0 & 0\\
c & d & 0 & 0 & 0 & 0 & 0\\
0 & 0 & a^{\prime\prime} & b^{\prime\prime} & 0 & 0 & 0\\
0 & 0 & c^{\prime\prime} & d^{\prime\prime} & 0 & 0 & 0\\
0 & 0 & 0 & 0 & e & 0 & 0\\
0 & 0 & 0 & 0 & 0 & d^{\prime} & c^{\prime}\\
0 & 0 & 0 & 0 & 0 & b^{\prime} & a^{\prime}%
\end{array}
\right)  . \label{sol.det.anotherpattern.a.Agg=}%
\end{equation}
But Corollary \ref{cor.sol.det.anotherpattern.gammagamma} (applied to $n=7$)
yields $\det\left(  A_{\left[  \gamma,\gamma\right]  }\right)  =\det A$.
Hence,%
\begin{align}
\det A  &  =\det\left(  A_{\left[  \gamma,\gamma\right]  }\right)
=\det\left(
\begin{array}
[c]{ccccccc}%
a & b & 0 & 0 & 0 & 0 & 0\\
c & d & 0 & 0 & 0 & 0 & 0\\
0 & 0 & a^{\prime\prime} & b^{\prime\prime} & 0 & 0 & 0\\
0 & 0 & c^{\prime\prime} & d^{\prime\prime} & 0 & 0 & 0\\
0 & 0 & 0 & 0 & e & 0 & 0\\
0 & 0 & 0 & 0 & 0 & d^{\prime} & c^{\prime}\\
0 & 0 & 0 & 0 & 0 & b^{\prime} & a^{\prime}%
\end{array}
\right)  \ \ \ \ \ \ \ \ \ \ \left(  \text{by
(\ref{sol.det.anotherpattern.a.Agg=})}\right) \nonumber\\
&  =\det\left(
\begin{array}
[c]{cccc}%
a & b & 0 & 0\\
c & d & 0 & 0\\
0 & 0 & a^{\prime\prime} & b^{\prime\prime}\\
0 & 0 & c^{\prime\prime} & d^{\prime\prime}%
\end{array}
\right)  \cdot\det\left(
\begin{array}
[c]{ccc}%
e & 0 & 0\\
0 & d^{\prime} & c^{\prime}\\
0 & b^{\prime} & a^{\prime}%
\end{array}
\right)  \label{sol.det.anotherpattern.a.2}%
\end{align}
(by Exercise \ref{exe.block2x2.tridet}, applied to $\left(
\begin{array}
[c]{cccc}%
a & b & 0 & 0\\
c & d & 0 & 0\\
0 & 0 & a^{\prime\prime} & b^{\prime\prime}\\
0 & 0 & c^{\prime\prime} & d^{\prime\prime}%
\end{array}
\right)  $, $\left(
\begin{array}
[c]{ccc}%
0 & 0 & 0\\
0 & 0 & 0\\
0 & 0 & 0\\
0 & 0 & 0
\end{array}
\right)  $, $\left(
\begin{array}
[c]{ccc}%
e & 0 & 0\\
0 & d^{\prime} & c^{\prime}\\
0 & b^{\prime} & a^{\prime}%
\end{array}
\right)  $, $3$ and $4$ instead of $A$, $B$, $D$, $m$ and $n$).

But Exercise \ref{exe.block2x2.tridet} (applied to $\left(
\begin{array}
[c]{cc}%
a & b\\
c & d
\end{array}
\right)  $, $\left(
\begin{array}
[c]{cc}%
0 & 0\\
0 & 0
\end{array}
\right)  $, $\left(
\begin{array}
[c]{cc}%
a^{\prime\prime} & b^{\prime\prime}\\
c^{\prime\prime} & d^{\prime\prime}%
\end{array}
\right)  $, $2$ and $2$ instead of $A$, $B$, $D$, $m$ and $n$) yields%
\[
\det\left(
\begin{array}
[c]{cccc}%
a & b & 0 & 0\\
c & d & 0 & 0\\
0 & 0 & a^{\prime\prime} & b^{\prime\prime}\\
0 & 0 & c^{\prime\prime} & d^{\prime\prime}%
\end{array}
\right)  =\underbrace{\det\left(
\begin{array}
[c]{cc}%
a & b\\
c & d
\end{array}
\right)  }_{=ad-bc}\cdot\underbrace{\det\left(
\begin{array}
[c]{cc}%
a^{\prime\prime} & b^{\prime\prime}\\
c^{\prime\prime} & d^{\prime\prime}%
\end{array}
\right)  }_{=a^{\prime\prime}d^{\prime\prime}-b^{\prime\prime}c^{\prime\prime
}}=\left(  ad-bc\right)  \cdot\left(  a^{\prime\prime}d^{\prime\prime
}-b^{\prime\prime}c^{\prime\prime}\right)  .
\]
Also, Exercise \ref{exe.block2x2.tridet} (applied to $\left(
\begin{array}
[c]{c}%
e
\end{array}
\right)  $, $\left(
\begin{array}
[c]{cc}%
0 & 0
\end{array}
\right)  $, $\left(
\begin{array}
[c]{cc}%
d^{\prime} & c^{\prime}\\
b^{\prime} & a^{\prime}%
\end{array}
\right)  $, $2$ and $1$ instead of $A$, $B$, $D$, $m$ and $n$) yields%
\[
\det\left(
\begin{array}
[c]{ccc}%
e & 0 & 0\\
0 & d^{\prime} & c^{\prime}\\
0 & b^{\prime} & a^{\prime}%
\end{array}
\right)  =\underbrace{\det\left(
\begin{array}
[c]{c}%
e
\end{array}
\right)  }_{=e}\cdot\underbrace{\det\left(
\begin{array}
[c]{cc}%
d^{\prime} & c^{\prime}\\
b^{\prime} & a^{\prime}%
\end{array}
\right)  }_{\substack{=d^{\prime}a^{\prime}-c^{\prime}b^{\prime}\\=a^{\prime
}d^{\prime}-b^{\prime}c^{\prime}}}=e\cdot\left(  a^{\prime}d^{\prime
}-b^{\prime}c^{\prime}\right)  .
\]
Hence, (\ref{sol.det.anotherpattern.a.2}) becomes%
\begin{align*}
\det A  &  =\underbrace{\det\left(
\begin{array}
[c]{cccc}%
a & b & 0 & 0\\
c & d & 0 & 0\\
0 & 0 & a^{\prime\prime} & b^{\prime\prime}\\
0 & 0 & c^{\prime\prime} & d^{\prime\prime}%
\end{array}
\right)  }_{=\left(  ad-bc\right)  \cdot\left(  a^{\prime\prime}%
d^{\prime\prime}-b^{\prime\prime}c^{\prime\prime}\right)  }\cdot
\underbrace{\det\left(
\begin{array}
[c]{ccc}%
e & 0 & 0\\
0 & d^{\prime} & c^{\prime}\\
0 & b^{\prime} & a^{\prime}%
\end{array}
\right)  }_{=e\cdot\left(  a^{\prime}d^{\prime}-b^{\prime}c^{\prime}\right)
}\\
&  =\left(  ad-bc\right)  \cdot\left(  a^{\prime\prime}d^{\prime\prime
}-b^{\prime\prime}c^{\prime\prime}\right)  \cdot e\cdot\left(  a^{\prime
}d^{\prime}-b^{\prime}c^{\prime}\right) \\
&  =e\cdot\left(  ad-bc\right)  \cdot\left(  a^{\prime\prime}d^{\prime\prime
}-b^{\prime\prime}c^{\prime\prime}\right)  \cdot\left(  a^{\prime}d^{\prime
}-b^{\prime}c^{\prime}\right)  .
\end{align*}
This solves Exercise \ref{exe.det.anotherpattern} \textbf{(a)}.

\textbf{(b)} Let $a,b,c,d,e,f,g,h,k,\ell,m,n$ be elements of $\mathbb{K}$. Let
$A$ be the $6\times6$-matrix%
\[
\left(
\begin{array}
[c]{cccccc}%
a & 0 & 0 & \ell & 0 & 0\\
0 & b & 0 & 0 & m & 0\\
0 & 0 & c & 0 & 0 & n\\
g & 0 & 0 & d & 0 & 0\\
0 & h & 0 & 0 & e & 0\\
0 & 0 & k & 0 & 0 & f
\end{array}
\right)  .
\]
Thus, the exercise demands that we find $\det A$.

Let $\gamma\in S_{6}$ be the permutation of the set $\left\{
1,2,3,4,5,6\right\}  $ which sends $1,2,3,4,5,6$ to $1,4,2,5,3,6$,
respectively. Then, straightforward computation (using the definitions of $A$
and of $A_{\left[  \gamma,\gamma\right]  }$) shows that%
\begin{equation}
A_{\left[  \gamma,\gamma\right]  }=\left(
\begin{array}
[c]{cccccc}%
a & \ell & 0 & 0 & 0 & 0\\
g & d & 0 & 0 & 0 & 0\\
0 & 0 & b & m & 0 & 0\\
0 & 0 & h & e & 0 & 0\\
0 & 0 & 0 & 0 & c & n\\
0 & 0 & 0 & 0 & k & f
\end{array}
\right)  . \label{sol.det.anotherpattern.b.Agg=}%
\end{equation}
But Corollary \ref{cor.sol.det.anotherpattern.gammagamma} (applied to $n=6$)
yields $\det\left(  A_{\left[  \gamma,\gamma\right]  }\right)  =\det A$.
Hence,%
\begin{align}
\det A  &  =\det\left(  A_{\left[  \gamma,\gamma\right]  }\right)
=\det\left(
\begin{array}
[c]{cccccc}%
a & \ell & 0 & 0 & 0 & 0\\
g & d & 0 & 0 & 0 & 0\\
0 & 0 & b & m & 0 & 0\\
0 & 0 & h & e & 0 & 0\\
0 & 0 & 0 & 0 & c & n\\
0 & 0 & 0 & 0 & k & f
\end{array}
\right)  \ \ \ \ \ \ \ \ \ \ \left(  \text{by
(\ref{sol.det.anotherpattern.b.Agg=})}\right) \nonumber\\
&  =\det\left(
\begin{array}
[c]{cccc}%
a & \ell & 0 & 0\\
g & d & 0 & 0\\
0 & 0 & b & m\\
0 & 0 & h & e
\end{array}
\right)  \cdot\det\left(
\begin{array}
[c]{cc}%
c & n\\
k & f
\end{array}
\right)  \label{sol.det.anotherpattern.b.1}%
\end{align}
(by Exercise \ref{exe.block2x2.tridet}, applied to $\left(
\begin{array}
[c]{cccc}%
a & \ell & 0 & 0\\
g & d & 0 & 0\\
0 & 0 & b & m\\
0 & 0 & h & e
\end{array}
\right)  $, $\left(
\begin{array}
[c]{cc}%
0 & 0\\
0 & 0\\
0 & 0\\
0 & 0
\end{array}
\right)  $, $\left(
\begin{array}
[c]{cc}%
c & n\\
k & f
\end{array}
\right)  $, $2$ and $4$ instead of $A$, $B$, $D$, $m$ and $n$).

But Exercise \ref{exe.block2x2.tridet} (applied to $\left(
\begin{array}
[c]{cc}%
a & \ell\\
g & d
\end{array}
\right)  $, $\left(
\begin{array}
[c]{cc}%
0 & 0\\
0 & 0
\end{array}
\right)  $, $\left(
\begin{array}
[c]{cc}%
b & m\\
h & e
\end{array}
\right)  $, $2$ and $2$ instead of $A$, $B$, $D$, $m$ and $n$) yields%
\[
\det\left(
\begin{array}
[c]{cccc}%
a & \ell & 0 & 0\\
g & d & 0 & 0\\
0 & 0 & b & m\\
0 & 0 & h & e
\end{array}
\right)  =\underbrace{\det\left(
\begin{array}
[c]{cc}%
a & \ell\\
g & d
\end{array}
\right)  }_{=ad-\ell g}\cdot\underbrace{\det\left(
\begin{array}
[c]{cc}%
b & m\\
h & e
\end{array}
\right)  }_{=be-mh}=\left(  ad-\ell g\right)  \cdot\left(  be-mh\right)  .
\]
Hence, (\ref{sol.det.anotherpattern.b.1}) becomes%
\[
\det A=\underbrace{\det\left(
\begin{array}
[c]{cccc}%
a & \ell & 0 & 0\\
g & d & 0 & 0\\
0 & 0 & b & m\\
0 & 0 & h & e
\end{array}
\right)  }_{=\left(  ad-\ell g\right)  \cdot\left(  be-mh\right)  }%
\cdot\underbrace{\det\left(
\begin{array}
[c]{cc}%
c & n\\
k & f
\end{array}
\right)  }_{=cf-nk}=\left(  ad-\ell g\right)  \cdot\left(  be-mh\right)
\cdot\left(  cf-nk\right)  .
\]
This solves Exercise \ref{exe.det.anotherpattern} \textbf{(b)}.
\end{proof}

\subsection{Solution to Exercise \ref{exe.adj(AB)}}

\begin{vershort}
Before we start solving this exercise, let us show some lemmas. The first of
them is a (somewhat disguised) particular case of the Cauchy-Binet formula:
\end{vershort}

\begin{verlong}
Before we start solving this exercise, let us show some lemmas. The first of
them formalizes the (intuitively obvious) fact that the elements of a finite
set of integers can be listed in increasing order in exactly one way:

\begin{lemma}
\label{lem.adj(AB).set.increase}Let $n\in\mathbb{N}$. Let $a_{1},a_{2}%
,\ldots,a_{n}$ be $n$ integers such that $a_{1}<a_{2}<\cdots<a_{n}$. Let
$b_{1},b_{2},\ldots,b_{n}$ be $n$ integers such that $b_{1}<b_{2}<\cdots
<b_{n}$. Assume that $\left\{  a_{1},a_{2},\ldots,a_{n}\right\}  =\left\{
b_{1},b_{2},\ldots,b_{n}\right\}  $. Then, $\left(  a_{1},a_{2},\ldots
,a_{n}\right)  =\left(  b_{1},b_{2},\ldots,b_{n}\right)  $.
\end{lemma}

Lemma \ref{lem.adj(AB).set.increase} can be easily derived from Theorem
\ref{thm.ind.inclist.unex} (by observing that $\left(  a_{1},a_{2}%
,\ldots,a_{n}\right)  $ and $\left(  b_{1},b_{2},\ldots,b_{n}\right)  $ are
increasing lists of the same set $\left\{  a_{1},a_{2},\ldots,a_{n}\right\}
=\left\{  b_{1},b_{2},\ldots,b_{n}\right\}  $). Let us nevertheless give a
different proof of Lemma \ref{lem.adj(AB).set.increase}:

\begin{proof}
[Proof of Lemma \ref{lem.adj(AB).set.increase}.]We have $b_{1}<b_{2}%
<\cdots<b_{n}$. In other words, if $p$ and $q$ are two elements of $\left\{
1,2,\ldots,n\right\}  $ satisfying $p<q$, then%
\begin{equation}
b_{p}<b_{q}. \label{pf.lem.adj(AB).set.increase.binc}%
\end{equation}
Hence, if $p$ and $q$ are two elements of $\left\{  1,2,\ldots,n\right\}  $
satisfying $b_{p}=b_{q}$, then%
\begin{equation}
p=q \label{pf.lem.adj(AB).set.increase.bdj}%
\end{equation}
\footnote{\textit{Proof of (\ref{pf.lem.adj(AB).set.increase.bdj}):} Let $p$
and $q$ be two elements of $\left\{  1,2,\ldots,n\right\}  $ satisfying
$b_{p}=b_{q}$. If we had $p<q$, then we would have $b_{p}<b_{q}$ (by
(\ref{pf.lem.adj(AB).set.increase.binc})), which would contradict $b_{p}%
=b_{q}$. Thus, we cannot have $p<q$. In other words, we must have $p\geq q$.
If we had $q<p$, then we would have $b_{q}<b_{p}$ (by
(\ref{pf.lem.adj(AB).set.increase.binc}) (applied to $q$ and $p$ instead of
$p$ and $q$)), which would contradict $b_{q}=b_{p}$. Thus, we cannot have
$q<p$. In other words, we must have $q\geq p$. Combined with $p\geq q$, this
shows that $p=q$. This proves (\ref{pf.lem.adj(AB).set.increase.bdj}).}.

We define a map $A:\left\{  1,2,\ldots,n\right\}  \rightarrow\left\{
1,2,\ldots,n\right\}  $ as follows:

Let $i\in\left\{  1,2,\ldots,n\right\}  $. Then, $a_{i}\in\left\{  a_{1}%
,a_{2},\ldots,a_{n}\right\}  =\left\{  b_{1},b_{2},\ldots,b_{n}\right\}  $.
Hence, there exists a $j\in\left\{  1,2,\ldots,n\right\}  $ satisfying
$a_{i}=b_{j}$. Moreover, this $j$ is unique\footnote{\textit{Proof.} Let
$j_{1}$ and $j_{2}$ be two elements $j\in\left\{  1,2,\ldots,n\right\}  $
satisfying $a_{i}=b_{j}$. We shall show that $j_{1}=j_{2}$.
\par
We know that $j_{1}$ is an element $j\in\left\{  1,2,\ldots,n\right\}  $
satisfying $a_{i}=b_{j}$. In other words, $j_{1}$ is an element of $\left\{
1,2,\ldots,n\right\}  $ satisfying $a_{i}=b_{j_{1}}$.
\par
We know that $j_{2}$ is an element $j\in\left\{  1,2,\ldots,n\right\}  $
satisfying $a_{i}=b_{j}$. In other words, $j_{2}$ is an element of $\left\{
1,2,\ldots,n\right\}  $ satisfying $a_{i}=b_{j_{2}}$.
\par
We have $b_{j_{1}}=a_{i}=b_{j_{2}}$. Hence, $j_{1}=j_{2}$ (by
(\ref{pf.lem.adj(AB).set.increase.bdj}), applied to $p=j_{1}$ and $q=j_{2}$).
\par
Now, let us forget that we fixed $j_{1}$ and $j_{2}$. We thus have proven that
if $j_{1}$ and $j_{2}$ are two elements $j\in\left\{  1,2,\ldots,n\right\}  $
satisfying $a_{i}=b_{j}$, then $j_{1}=j_{2}$. In other words, any two elements
$j\in\left\{  1,2,\ldots,n\right\}  $ satisfying $a_{i}=b_{j}$ must be equal.
Hence, the $j\in\left\{  1,2,\ldots,n\right\}  $ satisfying $a_{i}=b_{j}$ is
unique (because we already know that such a $j$ exists). Qed.}. We define
$A\left(  i\right)  $ to be this $j$. Thus, $A\left(  i\right)  $ is the
unique $j\in\left\{  1,2,\ldots,n\right\}  $ satisfying $a_{i}=b_{j}$. In
other words, $A\left(  i\right)  $ is an element of $\left\{  1,2,\ldots
,n\right\}  $ satisfying $a_{i}=b_{A\left(  i\right)  }$.

Thus, we have defined $A\left(  i\right)  \in\left\{  1,2,\ldots,n\right\}  $
for every $i\in\left\{  1,2,\ldots,n\right\}  $. In other words, we have
defined a map $A:\left\{  1,2,\ldots,n\right\}  \rightarrow\left\{
1,2,\ldots,n\right\}  $.

Thus, we have defined a map $A:\left\{  1,2,\ldots,n\right\}  \rightarrow
\left\{  1,2,\ldots,n\right\}  $ which satisfies%
\begin{equation}
\left(  a_{i}=b_{A\left(  i\right)  }\ \ \ \ \ \ \ \ \ \ \text{for every }%
i\in\left\{  1,2,\ldots,n\right\}  \right)  .
\label{pf.lem.adj(AB).set.increase.A}%
\end{equation}
An analogous construction (with $\left(  a_{1},a_{2},\ldots,a_{n}\right)  $,
$\left(  b_{1},b_{2},\ldots,b_{n}\right)  $ and $A$ replaced by $\left(
b_{1},b_{2},\ldots,b_{n}\right)  $, $\left(  a_{1},a_{2},\ldots,a_{n}\right)
$ and $B$) allows us to construct a map $B:\left\{  1,2,\ldots,n\right\}
\rightarrow\left\{  1,2,\ldots,n\right\}  $ which satisfies%
\begin{equation}
\left(  b_{i}=a_{B\left(  i\right)  }\ \ \ \ \ \ \ \ \ \ \text{for every }%
i\in\left\{  1,2,\ldots,n\right\}  \right)  .
\label{pf.lem.adj(AB).set.increase.B}%
\end{equation}
Consider this map $B$.

We have $A\circ B=\operatorname*{id}$\ \ \ \ \footnote{\textit{Proof.} Let
$i\in\left\{  1,2,\ldots,n\right\}  $. Then, $a_{B\left(  i\right)
}=b_{A\left(  B\left(  i\right)  \right)  }$ (by
(\ref{pf.lem.adj(AB).set.increase.A}), applied to $B\left(  i\right)  $
instead of $i$). Hence, $b_{A\left(  B\left(  i\right)  \right)  }=a_{B\left(
i\right)  }=b_{i}$ (by (\ref{pf.lem.adj(AB).set.increase.B})). Thus, $A\left(
B\left(  i\right)  \right)  =i$ (by (\ref{pf.lem.adj(AB).set.increase.bdj}),
applied to $p=A\left(  B\left(  i\right)  \right)  $ and $q=i$). Thus,
$\left(  A\circ B\right)  \left(  i\right)  =A\left(  B\left(  i\right)
\right)  =i=\operatorname*{id}\left(  i\right)  $.
\par
Now, let us forget that we fixed $i$. Thus, we have shown that $\left(  A\circ
B\right)  \left(  i\right)  =\operatorname*{id}\left(  i\right)  $ for every
$i\in\left\{  1,2,\ldots,n\right\}  $. In other words, $A\circ
B=\operatorname*{id}$, qed.}. The same argument (but applied to $\left(
b_{1},b_{2},\ldots,b_{n}\right)  $, $\left(  a_{1},a_{2},\ldots,a_{n}\right)
$, $B$ and $A$ instead of $\left(  a_{1},a_{2},\ldots,a_{n}\right)  $,
$\left(  b_{1},b_{2},\ldots,b_{n}\right)  $, $A$ and $B$) shows that $B\circ
A=\operatorname*{id}$.

The maps $A$ and $B$ are mutually inverse (since $A\circ B=\operatorname*{id}$
and $B\circ A=\operatorname*{id}$). Thus, the map $A$ is invertible. In other
words, the map $A$ is a bijection. Hence, $A$ is a bijection $\left\{
1,2,\ldots,n\right\}  \rightarrow\left\{  1,2,\ldots,n\right\}  $. In other
words, $A$ is a permutation of the set $\left\{  1,2,\ldots,n\right\}  $. In
other words, $A\in S_{n}$ (since $S_{n}$ is the set of all permutations of the
set $\left\{  1,2,\ldots,n\right\}  $).

The integers $b_{1},b_{2},\ldots,b_{n}$ are distinct (because of
(\ref{pf.lem.adj(AB).set.increase.bdj})). Proposition \ref{prop.sorting}
\textbf{(c)} (applied to $\left(  b_{1},b_{2},\ldots,b_{n}\right)  $ instead
of $\left(  a_{1},a_{2},\ldots,a_{n}\right)  $) thus shows that there is a
\textbf{unique} permutation $\sigma\in S_{n}$ such that $b_{\sigma\left(
1\right)  }<b_{\sigma\left(  2\right)  }<\cdots<b_{\sigma\left(  n\right)  }$.
Hence, there exists at most one such permutation. In other words, if
$\sigma_{1}$ and $\sigma_{2}$ are two permutations $\sigma\in S_{n}$ such that
$b_{\sigma\left(  1\right)  }<b_{\sigma\left(  2\right)  }<\cdots
<b_{\sigma\left(  n\right)  }$, then%
\begin{equation}
\sigma_{1}=\sigma_{2}. \label{pf.lem.adj(AB).set.increase.uni}%
\end{equation}

The permutation $\operatorname*{id}\in S_{n}$ satisfies $b_{\operatorname*{id}%
\left(  1\right)  }<b_{\operatorname*{id}\left(  2\right)  }<\cdots
<b_{\operatorname*{id}\left(  n\right)  }$\ \ \ \ \footnote{In fact, this is
just a way to rewrite the chain of inequalities $b_{1}<b_{2}<\cdots<b_{n}$
(which is true by assumption).}. In other words, the permutation
$\operatorname*{id}\in S_{n}$ is a permutation $\sigma\in S_{n}$ such that
$b_{\sigma\left(  1\right)  }<b_{\sigma\left(  2\right)  }<\cdots
<b_{\sigma\left(  n\right)  }$.

We have $a_{1}<a_{2}<\cdots<a_{n}$. This rewrites as $b_{A\left(  1\right)
}<b_{A\left(  2\right)  }<\cdots<b_{A\left(  n\right)  }$ (since
$a_{i}=b_{A\left(  i\right)  }$ for every $i\in\left\{  1,2,\ldots,n\right\}
$ (by (\ref{pf.lem.adj(AB).set.increase.A}))). Thus, $A$ is a permutation
$\sigma\in S_{n}$ such that $b_{\sigma\left(  1\right)  }<b_{\sigma\left(
2\right)  }<\cdots<b_{\sigma\left(  n\right)  }$.

Thus, we know that $A$ and $\operatorname*{id}$ are two permutations
$\sigma\in S_{n}$ such that $b_{\sigma\left(  1\right)  }<b_{\sigma\left(
2\right)  }<\cdots<b_{\sigma\left(  n\right)  }$. Hence,
(\ref{pf.lem.adj(AB).set.increase.uni}) (applied to $\sigma_{1}=A$ and
$\sigma_{2}=\operatorname*{id}$) shows that $A=\operatorname*{id}$.

But (\ref{pf.lem.adj(AB).set.increase.A}) shows that%
\begin{align*}
\left(  a_{1},a_{2},\ldots,a_{n}\right)   &  =\left(  b_{A\left(  1\right)
},b_{A\left(  2\right)  },\ldots,b_{A\left(  n\right)  }\right) \\
&  =\left(  b_{\operatorname*{id}\left(  1\right)  },b_{\operatorname*{id}%
\left(  2\right)  },\ldots,b_{\operatorname*{id}\left(  n\right)  }\right)
\ \ \ \ \ \ \ \ \ \ \left(  \text{since }A=\operatorname*{id}\right) \\
&  =\left(  b_{1},b_{2},\ldots,b_{n}\right)  \ \ \ \ \ \ \ \ \ \ \left(
\text{since }b_{\operatorname*{id}\left(  i\right)  }=b_{i}\text{ for every
}i\in\left\{  1,2,\ldots,n\right\}  \right)  .
\end{align*}
This proves Lemma \ref{lem.adj(AB).set.increase}.
\end{proof}
\end{verlong}

\begin{lemma}
\label{lem.adj(AB).cauchy-binet}Let $n$ be a positive integer. Let $A$ be an
$\left(  n-1\right)  \times n$-matrix. Let $B$ be an $n\times\left(
n-1\right)  $-matrix. Then,%
\[
\det\left(  AB\right)  =\sum_{k=1}^{n}\det\left(  \operatorname*{cols}%
\nolimits_{1,2,\ldots,\widehat{k},\ldots,n}A\right)  \cdot\det\left(
\operatorname*{rows}\nolimits_{1,2,\ldots,\widehat{k},\ldots,n}B\right)  .
\]

\end{lemma}

\begin{vershort}
\begin{proof}
[Proof of Lemma \ref{lem.adj(AB).cauchy-binet}.]Let $\left[  n\right]  $
denote the set $\left\{  1,2,\ldots,n\right\}  $. Define a subset $\mathbf{I}$
of $\left[  n\right]  ^{n-1}$ by%
\[
\mathbf{I}=\left\{  \left(  k_{1},k_{2},\ldots,k_{n-1}\right)  \in\left[
n\right]  ^{n-1}\ \mid\ k_{1}<k_{2}<\cdots<k_{n-1}\right\}  .
\]

Theorem \ref{thm.cauchy-binet} (applied to $n-1$ and $n$ instead of $n$ and
$m$) yields%
\begin{align}
&  \det\left(  AB\right) \nonumber\\
&  =\sum_{1\leq g_{1}<g_{2}<\cdots<g_{n-1}\leq n}\det\left(
\operatorname*{cols}\nolimits_{g_{1},g_{2},\ldots,g_{n-1}}A\right)  \cdot
\det\left(  \operatorname*{rows}\nolimits_{g_{1},g_{2},\ldots,g_{n-1}%
}B\right)  . \label{pf.lem.adj(AB).cauchy-binet.short.1}%
\end{align}
Recall that the summation sign $\sum_{1\leq g_{1}<g_{2}<\cdots<g_{n-1}\leq
n-1}$ is an abbreviation for \newline$\sum_{\substack{\left(  g_{1}%
,g_{2},\ldots,g_{n-1}\right)  \in\left\{  1,2,\ldots,n\right\}  ^{n-1}%
;\\g_{1}<g_{2}<\cdots<g_{n-1}}}$, which can be rewritten as $\sum_{\left(
g_{1},g_{2},\ldots,g_{n-1}\right)  \in\mathbf{I}}$ (because the $\left(
n-1\right)  $-tuples $\left(  g_{1},g_{2},\ldots,g_{n-1}\right)  \in\left\{
1,2,\ldots,n\right\}  ^{n-1}$ satisfying $g_{1}<g_{2}<\cdots<g_{n-1}$ are
precisely the elements of $\mathbf{I}$). Therefore,
(\ref{pf.lem.adj(AB).cauchy-binet.short.1}) can be rewritten as%
\begin{align}
&  \det\left(  AB\right) \nonumber\\
&  =\sum_{\left(  g_{1},g_{2},\ldots,g_{n-1}\right)  \in\mathbf{I}}\det\left(
\operatorname*{cols}\nolimits_{g_{1},g_{2},\ldots,g_{n-1}}A\right)  \cdot
\det\left(  \operatorname*{rows}\nolimits_{g_{1},g_{2},\ldots,g_{n-1}%
}B\right)  . \label{pf.lem.adj(AB).cauchy-binet.short.2}%
\end{align}

Now, let us take a closer look at $\mathbf{I}$. The set $\mathbf{I}$ consists
of all $\left(  n-1\right)  $-tuples $\left(  k_{1},k_{2},\ldots
,k_{n-1}\right)  \in\left[  n\right]  ^{n-1}$ satisfying $k_{1}<k_{2}%
<\cdots<k_{n-1}$. There are only $n$ such $\left(  n-1\right)  $-tuples:
namely, the $\left(  n-1\right)  $-tuples $\left(  1,2,\ldots,\widehat{k}%
,\ldots,n\right)  $ for $k\in\left\{  1,2,\ldots,n\right\}  $. This is
intuitively clear: If you want to choose an $\left(  n-1\right)  $-tuple
$\left(  k_{1},k_{2},\ldots,k_{n-1}\right)  \in\mathbf{I}$, you can simply
decide which of the $n$ elements $1,2,\ldots,n$ you do \textbf{not} want to be
an entry of $\left(  k_{1},k_{2},\ldots,k_{n-1}\right)  $, and then the
$\left(  n-1\right)  $-tuple $\left(  k_{1},k_{2},\ldots,k_{n-1}\right)  $
will have to be the list of all the remaining $n-1$ elements of $\left\{
1,2,\ldots,n\right\}  $ in increasing order. Let us formalize this argument a
bit more:

For every $k\in\left\{  1,2,\ldots,n\right\}  $, we have $\left(
1,2,\ldots,\widehat{k},\ldots,n\right)  \in\mathbf{I}$ (for obvious reasons).
Hence, we can define a map%
\[
\Phi:\left\{  1,2,\ldots,n\right\}  \rightarrow\mathbf{I}%
\]
by%
\[
\left(  \Phi\left(  k\right)  =\left(  1,2,\ldots,\widehat{k},\ldots,n\right)
\ \ \ \ \ \ \ \ \ \ \text{for every }k\in\left\{  1,2,\ldots,n\right\}
\right)  .
\]
Consider this map $\Phi$. This map $\Phi$ is
injective\footnote{\textit{Proof.} Let $i\in\left\{  1,2,\ldots,n\right\}  $
and $j\in\left\{  1,2,\ldots,n\right\}  $ be such that $\Phi\left(  i\right)
=\Phi\left(  j\right)  $. We shall show that $i=j$.
\par
The definition of $\Phi$ yields $\Phi\left(  i\right)  =\left(  1,2,\ldots
,\widehat{i},\ldots,n\right)  $. Hence, $i$ is the only element of $\left\{
1,2,\ldots,n\right\}  $ that does not appear in $\Phi\left(  i\right)  $.
Similarly, $j$ is the only element of $\left\{  1,2,\ldots,n\right\}  $ that
does not appear in $\Phi\left(  j\right)  $. In other words, $j$ is the only
element of $\left\{  1,2,\ldots,n\right\}  $ that does not appear in
$\Phi\left(  i\right)  $ (since $\Phi\left(  i\right)  =\Phi\left(  j\right)
$). Comparing this with the fact that $i$ is the only element of $\left\{
1,2,\ldots,n\right\}  $ that does not appear in $\Phi\left(  i\right)  $, we
conclude that $i=j$.
\par
Now, let us forget that we fixed $i$ and $j$. We thus have proven that if
$i\in\left\{  1,2,\ldots,n\right\}  $ and $j\in\left\{  1,2,\ldots,n\right\}
$ are such that $\Phi\left(  i\right)  =\Phi\left(  j\right)  $, then $i=j$.
In other words, the map $\Phi$ is injective.} and
surjective\footnote{\textit{Proof.} Let $\mathbf{g}\in\mathbf{I}$. We shall
show that $\mathbf{g}\in\Phi\left(  \left\{  1,2,\ldots,n\right\}  \right)  $.
\par
We have $\mathbf{g}\in\mathbf{I}=\left\{  \left(  k_{1},k_{2},\ldots
,k_{n-1}\right)  \in\left[  n\right]  ^{n-1}\ \mid\ k_{1}<k_{2}<\cdots
<k_{n-1}\right\}  $. In other words, $\mathbf{g}$ can be written in the form
$\mathbf{g}=\left(  g_{1},g_{2},\ldots,g_{n-1}\right)  $ for some $\left(
g_{1},g_{2},\ldots,g_{n-1}\right)  \in\left[  n\right]  ^{n-1}$ satisfying
$g_{1}<g_{2}<\cdots<g_{n-1}$. Consider this $\left(  g_{1},g_{2}%
,\ldots,g_{n-1}\right)  $.
\par
The integers $g_{1},g_{2},\ldots,g_{n-1}$ are distinct (since $g_{1}%
<g_{2}<\cdots<g_{n-1}$). Thus, $\left\{  g_{1},g_{2},\ldots,g_{n-1}\right\}  $
is an $\left(  n-1\right)  $-element subset of $\left[  n\right]  $.
Therefore, its complement $\left[  n\right]  \setminus\left\{  g_{1}%
,g_{2},\ldots,g_{n-1}\right\}  $ is a $1$-element subset of $\left[  n\right]
$ (since $n-\left(  n-1\right)  =1$). In other words, $\left[  n\right]
\setminus\left\{  g_{1},g_{2},\ldots,g_{n-1}\right\}  =\left\{  k\right\}  $
for some $k\in\left[  n\right]  $. Consider this $k$.
\par
We have $k\in\left[  n\right]  =\left\{  1,2,\ldots,n\right\}  $. Since
$\left\{  g_{1},g_{2},\ldots,g_{n-1}\right\}  \subseteq\left[  n\right]  $, we
have%
\[
\left\{  g_{1},g_{2},\ldots,g_{n-1}\right\}  =\left[  n\right]  \setminus
\underbrace{\left(  \left[  n\right]  \setminus\left\{  g_{1},g_{2}%
,\ldots,g_{n-1}\right\}  \right)  }_{=\left\{  k\right\}  }=\left[  n\right]
\setminus\left\{  k\right\}  .
\]
\par
Now, recall that $g_{1}<g_{2}<\cdots<g_{n-1}$. Hence, the $\left(  n-1\right)
$-tuple $\left(  g_{1},g_{2},\ldots,g_{n-1}\right)  $ is the list of all
elements of the set $\left\{  g_{1},g_{2},\ldots,g_{n-1}\right\}  $ in
increasing order. Since $\left\{  g_{1},g_{2},\ldots,g_{n-1}\right\}  =\left[
n\right]  \setminus\left\{  k\right\}  $, this rewrites as follows: The
$\left(  n-1\right)  $-tuple $\left(  g_{1},g_{2},\ldots,g_{n-1}\right)  $ is
the list of all elements of the set $\left[  n\right]  \setminus\left\{
k\right\}  $ in increasing order. But clearly the latter list is $\left(
1,2,\ldots,\widehat{k},\ldots,n\right)  $. Thus, the $\left(  n-1\right)
$-tuple $\left(  g_{1},g_{2},\ldots,g_{n-1}\right)  $ is the list $\left(
1,2,\ldots,\widehat{k},\ldots,n\right)  $. In other words, $\left(
g_{1},g_{2},\ldots,g_{n-1}\right)  =\left(  1,2,\ldots,\widehat{k}%
,\ldots,n\right)  $, so that%
\[
\mathbf{g}=\left(  g_{1},g_{2},\ldots,g_{n-1}\right)  =\left(  1,2,\ldots
,\widehat{k},\ldots,n\right)  =\Phi\left(  k\right)  \in\Phi\left(  \left\{
1,2,\ldots,n\right\}  \right)  .
\]
\par
Now, let us forget that we fixed $\mathbf{g}$. We thus have proven that
$\mathbf{g}\in\Phi\left(  \left\{  1,2,\ldots,n\right\}  \right)  $ for every
$\mathbf{g}\in\mathbf{I}$. In other words, $\mathbf{I}\subseteq\Phi\left(
\left\{  1,2,\ldots,n\right\}  \right)  $. In other words, the map $\Phi$ is
surjective.}. Hence, the map $\Phi$ is a bijection. In other words, the map%
\[
\left\{  1,2,\ldots,n\right\}  \rightarrow\mathbf{I}%
,\ \ \ \ \ \ \ \ \ \ k\mapsto\left(  1,2,\ldots,\widehat{k},\ldots,n\right)
\]
is a bijection (since this map is precisely $\Phi$).

Now, (\ref{pf.lem.adj(AB).cauchy-binet.short.2}) becomes%
\begin{align*}
\det\left(  AB\right)   &  =\sum_{\left(  g_{1},g_{2},\ldots,g_{n-1}\right)
\in\mathbf{I}}\det\left(  \operatorname*{cols}\nolimits_{g_{1},g_{2}%
,\ldots,g_{n-1}}A\right)  \cdot\det\left(  \operatorname*{rows}%
\nolimits_{g_{1},g_{2},\ldots,g_{n-1}}B\right) \\
&  =\underbrace{\sum_{k\in\left\{  1,2,\ldots,n\right\}  }}_{=\sum_{k=1}^{n}%
}\det\left(  \operatorname*{cols}\nolimits_{1,2,\ldots,\widehat{k},\ldots
,n}A\right)  \cdot\det\left(  \operatorname*{rows}\nolimits_{1,2,\ldots
,\widehat{k},\ldots,n}B\right) \\
&  \ \ \ \ \ \ \ \ \ \ \left(
\begin{array}
[c]{c}%
\text{here, we have substituted }\left(  1,2,\ldots,\widehat{k},\ldots
,n\right)  \text{ for}\\
\left(  g_{1},g_{2},\ldots,g_{n-1}\right)  \text{ in the sum, since the map}\\
\left\{  1,2,\ldots,n\right\}  \rightarrow\mathbf{I},\ k\mapsto\left(
1,2,\ldots,\widehat{k},\ldots,n\right)  \text{ is a bijection}%
\end{array}
\right) \\
&  =\sum_{k=1}^{n}\det\left(  \operatorname*{cols}\nolimits_{1,2,\ldots
,\widehat{k},\ldots,n}A\right)  \cdot\det\left(  \operatorname*{rows}%
\nolimits_{1,2,\ldots,\widehat{k},\ldots,n}B\right)  .
\end{align*}
This proves Lemma \ref{lem.adj(AB).cauchy-binet}.
\end{proof}
\end{vershort}

\begin{verlong}
This lemma is a particular case of the Cauchy-Binet formula. Here is how it
can be derived from the latter, in detail:

\begin{proof}
[Proof of Lemma \ref{lem.adj(AB).cauchy-binet}.]Let $\left[  n\right]  $
denote the set $\left\{  1,2,\ldots,n\right\}  $. Define a subset $\mathbf{I}$
of $\left[  n\right]  ^{n-1}$ by%
\[
\mathbf{I}=\left\{  \left(  k_{1},k_{2},\ldots,k_{n-1}\right)  \in\left[
n\right]  ^{n-1}\ \mid\ k_{1}<k_{2}<\cdots<k_{n-1}\right\}  .
\]

Theorem \ref{thm.cauchy-binet} (applied to $n-1$ and $n$ instead of $n$ and
$m$) yields%
\begin{align}
&  \det\left(  AB\right) \nonumber\\
&  =\underbrace{\sum_{1\leq g_{1}<g_{2}<\cdots<g_{n-1}\leq n}}%
_{\substack{=\sum_{\substack{\left(  g_{1},g_{2},\ldots,g_{n-1}\right)
\in\left\{  1,2,\ldots,n\right\}  ^{n-1};\\g_{1}<g_{2}<\cdots<g_{n-1}}%
}\\=\sum_{\substack{\left(  g_{1},g_{2},\ldots,g_{n-1}\right)  \in\left[
n\right]  ^{n-1};\\g_{1}<g_{2}<\cdots<g_{n-1}}}\\\text{(since }\left\{
1,2,\ldots,n-1\right\}  =\left[  n\right]  \text{)}}}\det\left(
\operatorname*{cols}\nolimits_{g_{1},g_{2},\ldots,g_{n-1}}A\right)  \cdot
\det\left(  \operatorname*{rows}\nolimits_{g_{1},g_{2},\ldots,g_{n-1}}B\right)
\nonumber\\
&  =\underbrace{\sum_{\substack{\left(  g_{1},g_{2},\ldots,g_{n-1}\right)
\in\left[  n\right]  ^{n-1};\\g_{1}<g_{2}<\cdots<g_{n-1}}}}_{\substack{=\sum
_{\left(  g_{1},g_{2},\ldots,g_{n-1}\right)  \in\left\{  \left(  k_{1}%
,k_{2},\ldots,k_{n-1}\right)  \in\left[  n\right]  ^{n-1}\ \mid\ k_{1}%
<k_{2}<\cdots<k_{n-1}\right\}  }\\=\sum_{\left(  g_{1},g_{2},\ldots
,g_{n-1}\right)  \in\mathbf{I}}\\\text{(since }\left\{  \left(  k_{1}%
,k_{2},\ldots,k_{n-1}\right)  \in\left[  n\right]  ^{n-1}\ \mid\ k_{1}%
<k_{2}<\cdots<k_{n-1}\right\}  =\mathbf{I}\text{)}}}\det\left(
\operatorname*{cols}\nolimits_{g_{1},g_{2},\ldots,g_{n-1}}A\right)  \cdot
\det\left(  \operatorname*{rows}\nolimits_{g_{1},g_{2},\ldots,g_{n-1}}B\right)
\nonumber\\
&  =\sum_{\left(  g_{1},g_{2},\ldots,g_{n-1}\right)  \in\mathbf{I}}\det\left(
\operatorname*{cols}\nolimits_{g_{1},g_{2},\ldots,g_{n-1}}A\right)  \cdot
\det\left(  \operatorname*{rows}\nolimits_{g_{1},g_{2},\ldots,g_{n-1}%
}B\right)  . \label{pf.lem.adj(AB).cauchy-binet.1}%
\end{align}

Now, for every $k\in\left\{  1,2,\ldots,n\right\}  $, we have%
\[
\left(  1,2,\ldots,\widehat{k},\ldots,n\right)  \in\mathbf{I}%
\]
\footnote{\textit{Proof.} Let $k\in\left\{  1,2,\ldots,n\right\}  $. Let us
denote the $\left(  n-1\right)  $-tuple $\left(  1,2,\ldots,\widehat{k}%
,\ldots,n\right)  $ by $\left(  t_{1},t_{2},\ldots,t_{n-1}\right)  $. Thus,%
\begin{equation}
t_{i}=%
\begin{cases}
i, & \text{if }i<k;\\
i+1, & \text{if }i\geq k
\end{cases}
\ \ \ \ \ \ \ \ \ \ \text{for every }i\in\left\{  1,2,\ldots,n-1\right\}  .
\label{pf.lem.adj(AB).cauchy-binet.wd.pf.1}%
\end{equation}
\par
Now, let $i\in\left\{  1,2,\ldots,n-2\right\}  $. We shall prove that
$t_{i}<t_{i+1}$. Indeed, let us assume the contrary (for the sake of
contradiction). Thus, $t_{i}\geq t_{i+1}$. But%
\begin{align*}
t_{i}  &  =%
\begin{cases}
i, & \text{if }i<k;\\
i+1, & \text{if }i\geq k
\end{cases}
\ \ \ \ \ \ \ \ \ \ \left(  \text{by
(\ref{pf.lem.adj(AB).cauchy-binet.wd.pf.1})}\right) \\
&  \leq%
\begin{cases}
i+1, & \text{if }i<k;\\
i+1, & \text{if }i\geq k
\end{cases}
\ \ \ \ \ \ \ \ \ \ \left(  \text{since }i\leq i+1\text{ in the case when
}i<k\right) \\
&  =i+1,
\end{align*}
so that%
\begin{align*}
i+1  &  \geq t_{i}\geq t_{i+1}=%
\begin{cases}
i+1, & \text{if }i+1<k;\\
\left(  i+1\right)  +1, & \text{if }i+1\geq k
\end{cases}
\ \ \ \ \ \ \ \ \ \ \left(  \text{by
(\ref{pf.lem.adj(AB).cauchy-binet.wd.pf.1}), applied to }i+1\text{ instead of
}i\right) \\
&  \geq%
\begin{cases}
i+1, & \text{if }i+1<k;\\
i+1, & \text{if }i+1\geq k
\end{cases}
\ \ \ \ \ \ \ \ \ \ \left(  \text{since }\left(  i+1\right)  +1\geq i+1\text{
in the case when }i+1\geq k\right) \\
&  =i+1.
\end{align*}
Combining $i+1\geq t_{i}$ with $t_{i}\geq i+1$, we obtain $i+1=t_{i}$. Hence,
if we had $i<k$, then we would have%
\[
i+1=t_{i}=%
\begin{cases}
i, & \text{if }i<k;\\
i+1, & \text{if }i\geq k
\end{cases}
=i\ \ \ \ \ \ \ \ \ \ \left(  \text{since }i<k\right)  ,
\]
which would contradict $i+1\neq i$. Therefore, we cannot have $i<k$. Thus, we
have $i\geq k$. Thus, $i+1\geq k$. Now,%
\begin{align*}
t_{i+1}  &  =%
\begin{cases}
i+1, & \text{if }i+1<k;\\
\left(  i+1\right)  +1, & \text{if }i+1\geq k
\end{cases}
=\left(  i+1\right)  +1\ \ \ \ \ \ \ \ \ \ \left(  \text{since }i+1\geq
k\right) \\
&  >i+1=t_{i}\geq t_{i+1}.
\end{align*}
This is absurd. Thus, we have obtained a contradiction. This contradiction
proves that our assumption was wrong. Hence, $t_{i}<t_{i+1}$.
\par
Now, let us forget that we fixed $i$. We thus have shown that $t_{i}<t_{i+1}$
for every $i\in\left\{  1,2,\ldots,n-2\right\}  $. In other words,
$t_{1}<t_{2}<\cdots<t_{n-1}$.
\par
Also, $\left(  t_{1},t_{2},\ldots,t_{n-1}\right)  =\left(  1,2,\ldots
,\widehat{k},\ldots,n\right)  \in\left[  n\right]  ^{n-1}$. Hence, $\left(
t_{1},t_{2},\ldots,t_{n-1}\right)  $ is an element of $\left[  n\right]
^{n-1}$ and satisfies $t_{1}<t_{2}<\cdots<t_{n-1}$. In other words, $\left(
t_{1},t_{2},\ldots,t_{n-1}\right)  $ is an element $\left(  k_{1},k_{2}%
,\ldots,k_{n-1}\right)  \in\left[  n\right]  ^{n-1}$ satisfying $k_{1}%
<k_{2}<\cdots<k_{n-1}$. In other words,
\[
\left(  t_{1},t_{2},\ldots,t_{n-1}\right)  \in\left\{  \left(  k_{1}%
,k_{2},\ldots,k_{n-1}\right)  \in\left[  n\right]  ^{n-1}\ \mid\ k_{1}%
<k_{2}<\cdots<k_{n-1}\right\}  =\mathbf{I}.
\]
Thus, $\left(  1,2,\ldots,\widehat{k},\ldots,n\right)  =\left(  t_{1}%
,t_{2},\ldots,t_{n-1}\right)  \in\mathbf{I}$, qed.}. Hence, we can define a
map%
\[
\Phi:\left\{  1,2,\ldots,n\right\}  \rightarrow\mathbf{I}%
\]
by%
\[
\left(  \Phi\left(  k\right)  =\left(  1,2,\ldots,\widehat{k},\ldots,n\right)
\ \ \ \ \ \ \ \ \ \ \text{for every }k\in\left\{  1,2,\ldots,n\right\}
\right)  .
\]
Consider this map $\Phi$. This map $\Phi$ is
injective\footnote{\textit{Proof.} Let $i\in\left\{  1,2,\ldots,n\right\}  $
and $j\in\left\{  1,2,\ldots,n\right\}  $ be such that $\Phi\left(  i\right)
=\Phi\left(  j\right)  $. We shall show that $i=j$.
\par
The definition of $\Phi$ yields $\Phi\left(  i\right)  =\left(  1,2,\ldots
,\widehat{i},\ldots,n\right)  $. Hence,%
\begin{align*}
&  \left(  \text{the set of all entries of }\underbrace{\Phi\left(  i\right)
}_{=\left(  1,2,\ldots,\widehat{i},\ldots,n\right)  }\right) \\
&  =\left(  \text{the set of all entries of }\left(  1,2,\ldots,\widehat{i}%
,\ldots,n\right)  \right) \\
&  =\left\{  1,2,\ldots,\widehat{i},\ldots,n\right\}  =\left\{  1,2,\ldots
,n\right\}  \setminus\left\{  i\right\}  .
\end{align*}
Thus,%
\begin{align*}
&  \left\{  1,2,\ldots,n\right\}  \setminus\underbrace{\left(  \text{the set
of all entries of }\Phi\left(  i\right)  \right)  }_{=\left\{  1,2,\ldots
,n\right\}  \setminus\left\{  i\right\}  }\\
&  =\left\{  1,2,\ldots,n\right\}  \setminus\left(  \left\{  1,2,\ldots
,n\right\}  \setminus\left\{  i\right\}  \right)  =\left\{  i\right\}
\ \ \ \ \ \ \ \ \ \ \left(  \text{since }\left\{  i\right\}  \subseteq\left\{
1,2,\ldots,n\right\}  \text{ (since }i\in\left\{  1,2,\ldots,n\right\}
\text{)}\right)  .
\end{align*}
The same argument (applied to $j$ instead of $i$) shows that%
\[
\left\{  1,2,\ldots,n\right\}  \setminus\left(  \text{the set of all entries
of }\Phi\left(  j\right)  \right)  =\left\{  j\right\}  .
\]
\par
Now,%
\begin{align*}
\left\{  i\right\}   &  =\left\{  1,2,\ldots,n\right\}  \setminus\left(
\text{the set of all entries of }\underbrace{\Phi\left(  i\right)  }%
_{=\Phi\left(  j\right)  }\right) \\
&  =\left\{  1,2,\ldots,n\right\}  \setminus\left(  \text{the set of all
entries of }\Phi\left(  j\right)  \right)  =\left\{  j\right\}  .
\end{align*}
Hence, $i\in\left\{  i\right\}  =\left\{  j\right\}  $, so that $i=j$.
\par
Now, let us forget that we fixed $i$ and $j$. We thus have proven that if
$i\in\left\{  1,2,\ldots,n\right\}  $ and $j\in\left\{  1,2,\ldots,n\right\}
$ are such that $\Phi\left(  i\right)  =\Phi\left(  j\right)  $, then $i=j$.
In other words, the map $\Phi$ is injective, qed.} and
surjective\footnote{\textit{Proof.} Let $\mathbf{g}\in\mathbf{I}$. We shall
show that $\mathbf{g}\in\Phi\left(  \left\{  1,2,\ldots,n\right\}  \right)  $.
\par
We have%
\begin{align*}
\mathbf{g}  &  \in\mathbf{I}=\left\{  \left(  k_{1},k_{2},\ldots
,k_{n-1}\right)  \in\left[  n\right]  ^{n-1}\ \mid\ k_{1}<k_{2}<\cdots
<k_{n-1}\right\} \\
&  =\left\{  \left(  g_{1},g_{2},\ldots,g_{n-1}\right)  \in\left[  n\right]
^{n-1}\ \mid\ g_{1}<g_{2}<\cdots<g_{n-1}\right\}
\end{align*}
(here, we have renamed the index $\left(  k_{1},k_{2},\ldots,k_{n-1}\right)  $
as $\left(  g_{1},g_{2},\ldots,g_{n-1}\right)  $). In other words,
$\mathbf{g}$ can be written in the form $\mathbf{g}=\left(  g_{1},g_{2}%
,\ldots,g_{n-1}\right)  $ for some $\left(  g_{1},g_{2},\ldots,g_{n-1}\right)
\in\left[  n\right]  ^{n-1}$ satisfying $g_{1}<g_{2}<\cdots<g_{n-1}$. Consider
this $\left(  g_{1},g_{2},\ldots,g_{n-1}\right)  $.
\par
The integers $g_{1},g_{2},\ldots,g_{n-1}$ are distinct (since $g_{1}%
<g_{2}<\cdots<g_{n-1}$). Hence, they are $n-1$ distinct integers. In other
words, we have $\left\vert \left\{  g_{1},g_{2},\ldots,g_{n-1}\right\}
\right\vert =n-1$. Also, each $i\in\left\{  1,2,\ldots,n-1\right\}  $
satisfies $g_{i}\in\left[  n\right]  $ (since $\left(  g_{1},g_{2}%
,\ldots,g_{n-1}\right)  \in\left[  n\right]  ^{n-1}$). In other words,
$\left\{  g_{1},g_{2},\ldots,g_{n-1}\right\}  \subseteq\left[  n\right]  $.
Hence,%
\[
\left\vert \left[  n\right]  \setminus\left\{  g_{1},g_{2},\ldots
,g_{n-1}\right\}  \right\vert =\underbrace{\left\vert \left[  n\right]
\right\vert }_{=n}-\underbrace{\left\vert \left\{  g_{1},g_{2},\ldots
,g_{n-1}\right\}  \right\vert }_{=n-1}=n-\left(  n-1\right)  =1.
\]
In other words, $\left[  n\right]  \setminus\left\{  g_{1},g_{2}%
,\ldots,g_{n-1}\right\}  $ is a one-element set. Hence, the set $\left[
n\right]  \setminus\left\{  g_{1},g_{2},\ldots,g_{n-1}\right\}  $ has the form
$\left\{  k\right\}  $ for some object $k$. Consider this $k$.
\par
We have $\left[  n\right]  \setminus\left\{  g_{1},g_{2},\ldots,g_{n-1}%
\right\}  =\left\{  k\right\}  $, so that $k\in\left\{  k\right\}  =\left[
n\right]  \setminus\left\{  g_{1},g_{2},\ldots,g_{n-1}\right\}  \subseteq
\left[  n\right]  =\left\{  1,2,\ldots,n\right\}  $.
\par
Let us denote the $\left(  n-1\right)  $-tuple $\left(  1,2,\ldots
,\widehat{k},\ldots,n\right)  $ by $\left(  t_{1},t_{2},\ldots,t_{n-1}\right)
$. Thus,%
\begin{align*}
\left(  t_{1},t_{2},\ldots,t_{n-1}\right)   &  =\left(  1,2,\ldots
,\widehat{k},\ldots,n\right) \\
&  \in\mathbf{I}=\left\{  \left(  k_{1},k_{2},\ldots,k_{n-1}\right)
\in\left[  n\right]  ^{n-1}\ \mid\ k_{1}<k_{2}<\cdots<k_{n-1}\right\}  .
\end{align*}
In other words, $\left(  t_{1},t_{2},\ldots,t_{n-1}\right)  $ is an element
$\left(  k_{1},k_{2},\ldots,k_{n-1}\right)  $ of $\left[  n\right]  ^{n-1}$
satisfying $k_{1}<k_{2}<\cdots<k_{n-1}$. In other words, $\left(  t_{1}%
,t_{2},\ldots,t_{n-1}\right)  $ is an element of $\left[  n\right]  ^{n-1}$
and satisfies $t_{1}<t_{2}<\cdots<t_{n-1}$.
\par
Clearly,
\begin{align*}
\left\{  t_{1},t_{2},\ldots,t_{n-1}\right\}   &  =\left(  \text{the set of all
entries of the list }\underbrace{\left(  t_{1},t_{2},\ldots,t_{n-1}\right)
}_{=\left(  1,2,\ldots,\widehat{k},\ldots,n\right)  }\right) \\
&  =\left(  \text{the set of all entries of the list }\left(  1,2,\ldots
,\widehat{k},\ldots,n\right)  \right) \\
&  =\left\{  1,2,\ldots,\widehat{k},\ldots,n\right\}  =\underbrace{\left\{
1,2,\ldots,n\right\}  }_{=\left[  n\right]  }\setminus\underbrace{\left\{
k\right\}  }_{=\left[  n\right]  \setminus\left\{  g_{1},g_{2},\ldots
,g_{n-1}\right\}  }\\
&  =\left[  n\right]  \setminus\left(  \left[  n\right]  \setminus\left\{
g_{1},g_{2},\ldots,g_{n-1}\right\}  \right)  =\left\{  g_{1},g_{2}%
,\ldots,g_{n-1}\right\}
\end{align*}
(since $\left\{  g_{1},g_{2},\ldots,g_{n-1}\right\}  \subseteq\left[
n\right]  $).
\par
Now, Lemma \ref{lem.adj(AB).set.increase} (applied to $n-1$, $\left(
t_{1},t_{2},\ldots,t_{n-1}\right)  $ and $\left(  g_{1},g_{2},\ldots
,g_{n-1}\right)  $ instead of $n$, $\left(  a_{1},a_{2},\ldots,a_{n}\right)  $
and $\left(  b_{1},b_{2},\ldots,b_{n}\right)  $) shows that $\left(
t_{1},t_{2},\ldots,t_{n-1}\right)  =\left(  g_{1},g_{2},\ldots,g_{n-1}\right)
$. Compared with $\mathbf{g}=\left(  g_{1},g_{2},\ldots,g_{n-1}\right)  $,
this yields%
\[
\mathbf{g}=\left(  t_{1},t_{2},\ldots,t_{n-1}\right)  =\left(  1,2,\ldots
,\widehat{k},\ldots,n\right)  =\Phi\left(  k\right)
\]
(since $\Phi\left(  k\right)  =\left(  1,2,\ldots,\widehat{k},\ldots,n\right)
$ (by the definition of $\Phi\left(  k\right)  $)). Thus, $\mathbf{g}%
=\Phi\left(  \underbrace{k}_{\in\left\{  1,2,\ldots,n\right\}  }\right)
\in\Phi\left(  \left\{  1,2,\ldots,n\right\}  \right)  $.
\par
Now, let us forget that we fixed $\mathbf{g}$. We thus have proven that
$\mathbf{g}\in\Phi\left(  \left\{  1,2,\ldots,n\right\}  \right)  $ for every
$\mathbf{g}\in\mathbf{I}$. In other words, $\mathbf{I}\subseteq\Phi\left(
\left\{  1,2,\ldots,n\right\}  \right)  $. In other words, the map $\Phi$ is
surjective, qed.}. Hence, the map $\Phi$ is bijective. In other words, the map
$\Phi$ is a bijection. In other words, the map%
\begin{align*}
\left\{  1,2,\ldots,n\right\}   &  \rightarrow\mathbf{I},\\
k  &  \mapsto\left(  1,2,\ldots,\widehat{k},\ldots,n\right)
\end{align*}
is a bijection (because the map $\Phi$ is the map
\begin{align*}
\left\{  1,2,\ldots,n\right\}   &  \rightarrow\mathbf{I},\\
k  &  \mapsto\left(  1,2,\ldots,\widehat{k},\ldots,n\right)
\end{align*}
(since we have $\Phi\left(  k\right)  =\left(  1,2,\ldots,\widehat{k}%
,\ldots,n\right)  $ for every $k\in\left\{  1,2,\ldots,n\right\}  $)).

Now, (\ref{pf.lem.adj(AB).cauchy-binet.1}) becomes%
\begin{align*}
\det\left(  AB\right)   &  =\sum_{\left(  g_{1},g_{2},\ldots,g_{n-1}\right)
\in\mathbf{I}}\det\left(  \operatorname*{cols}\nolimits_{g_{1},g_{2}%
,\ldots,g_{n-1}}A\right)  \cdot\det\left(  \operatorname*{rows}%
\nolimits_{g_{1},g_{2},\ldots,g_{n-1}}B\right) \\
&  =\underbrace{\sum_{k\in\left\{  1,2,\ldots,n\right\}  }}_{=\sum_{k=1}^{n}%
}\det\left(  \operatorname*{cols}\nolimits_{1,2,\ldots,\widehat{k},\ldots
,n}A\right)  \cdot\det\left(  \operatorname*{rows}\nolimits_{1,2,\ldots
,\widehat{k},\ldots,n}B\right) \\
&  \ \ \ \ \ \ \ \ \ \ \left(
\begin{array}
[c]{c}%
\text{here, we have substituted }\left(  1,2,\ldots,\widehat{k},\ldots
,n\right)  \text{ for}\\
\left(  g_{1},g_{2},\ldots,g_{n-1}\right)  \text{ in the sum, since the map}\\
\left\{  1,2,\ldots,n\right\}  \rightarrow\mathbf{I},\ k\mapsto\left(
1,2,\ldots,\widehat{k},\ldots,n\right)  \text{ is a bijection}%
\end{array}
\right) \\
&  =\sum_{k=1}^{n}\det\left(  \operatorname*{cols}\nolimits_{1,2,\ldots
,\widehat{k},\ldots,n}A\right)  \cdot\det\left(  \operatorname*{rows}%
\nolimits_{1,2,\ldots,\widehat{k},\ldots,n}B\right)  .
\end{align*}
This proves Lemma \ref{lem.adj(AB).cauchy-binet}.
\end{proof}
\end{verlong}

Here comes one more simple lemma:

\begin{lemma}
\label{lem.adj(AB).minor-of-AB}Let $n\in\mathbb{N}$, $m\in\mathbb{N}$ and
$p\in\mathbb{N}$. Let $A$ be an $n\times p$-matrix. Let $B$ be a $p\times
m$-matrix. Let $i_{1},i_{2},\ldots,i_{u}$ be some elements of $\left\{
1,2,\ldots,n\right\}  $. Let $j_{1},j_{2},\ldots,j_{v}$ be some elements of
$\left\{  1,2,\ldots,m\right\}  $. Then,%
\[
\operatorname*{sub}\nolimits_{i_{1},i_{2},\ldots,i_{u}}^{j_{1},j_{2}%
,\ldots,j_{v}}\left(  AB\right)  =\left(  \operatorname*{rows}\nolimits_{i_{1}%
,i_{2},\ldots,i_{u}}A\right)  \cdot\left(  \operatorname*{cols}%
\nolimits_{j_{1},j_{2},\ldots,j_{v}}B\right)  .
\]

\end{lemma}

\begin{proof}
[Proof of Lemma \ref{lem.adj(AB).minor-of-AB}.]Let us write the $n\times
p$-matrix $A$ in the form $A=\left(  a_{i,j}\right)  _{1\leq i\leq n,\ 1\leq
j\leq p}$. Let us write the $p\times m$-matrix $B$ in the form $B=\left(
b_{i,j}\right)  _{1\leq i\leq p,\ 1\leq j\leq m}$.

The definition of the product of two matrices yields $AB=\left(  \sum
_{k=1}^{p}a_{i,k}b_{k,j}\right)  _{1\leq i\leq n,\ 1\leq j\leq m}$ (since
$A=\left(  a_{i,j}\right)  _{1\leq i\leq n,\ 1\leq j\leq p}$ and $B=\left(
b_{i,j}\right)  _{1\leq i\leq p,\ 1\leq j\leq m}$). Thus, the definition of
$\operatorname*{sub}\nolimits_{i_{1},i_{2},\ldots,i_{u}}^{j_{1},j_{2}%
,\ldots,j_{v}}\left(  AB\right)  $ yields%
\begin{equation}
\operatorname*{sub}\nolimits_{i_{1},i_{2},\ldots,i_{u}}^{j_{1},j_{2}%
,\ldots,j_{v}}\left(  AB\right)  =\left(  \sum_{k=1}^{p}a_{i_{x},k}b_{k,j_{y}%
}\right)  _{1\leq x\leq u,\ 1\leq y\leq v}.
\label{pf.lem.adj(AB).minor-of-AB.1}%
\end{equation}

\begin{vershort}
On the other hand, we have%
\[
\operatorname*{rows}\nolimits_{i_{1},i_{2},\ldots,i_{u}}A=\left(  a_{i_{x}%
,j}\right)  _{1\leq x\leq u,\ 1\leq j\leq p}%
\]
(by the definition of $\operatorname*{rows}\nolimits_{i_{1},i_{2},\ldots
,i_{u}}A$, since $A=\left(  a_{i,j}\right)  _{1\leq i\leq n,\ 1\leq j\leq p}$)
and%
\[
\operatorname*{cols}\nolimits_{j_{1},j_{2},\ldots,j_{v}}B=\left(  b_{i,j_{y}%
}\right)  _{1\leq i\leq p,\ 1\leq y\leq v}%
\]
(by the definition of $\operatorname*{cols}\nolimits_{j_{1},j_{2},\ldots
,j_{v}}B$, since $B=\left(  b_{i,j}\right)  _{1\leq i\leq p,\ 1\leq j\leq m}%
$). The definition of the product of two matrices thus yields%
\[
\left(  \operatorname*{rows}\nolimits_{i_{1},i_{2},\ldots,i_{u}}A\right)
\cdot\left(  \operatorname*{cols}\nolimits_{j_{1},j_{2},\ldots,j_{v}}B\right)
=\left(  \sum_{k=1}^{p}a_{i_{x},k}b_{k,j_{y}}\right)  _{1\leq x\leq u,\ 1\leq
y\leq v}.
\]
Comparing this with (\ref{pf.lem.adj(AB).minor-of-AB.1}), we obtain%
\[
\operatorname*{sub}\nolimits_{i_{1},i_{2},\ldots,i_{u}}^{j_{1},j_{2}%
,\ldots,j_{v}}\left(  AB\right)  =\left(  \operatorname*{rows}\nolimits_{i_{1}%
,i_{2},\ldots,i_{u}}A\right)  \cdot\left(  \operatorname*{cols}%
\nolimits_{j_{1},j_{2},\ldots,j_{v}}B\right)  .
\]

\end{vershort}

\begin{verlong}
On the other hand, $A=\left(  a_{i,j}\right)  _{1\leq i\leq n,\ 1\leq j\leq
p}$. Hence, the definition of $\operatorname*{rows}\nolimits_{i_{1}%
,i_{2},\ldots,i_{u}}A$ yields%
\[
\operatorname*{rows}\nolimits_{i_{1},i_{2},\ldots,i_{u}}A=\left(  a_{i_{x}%
,j}\right)  _{1\leq x\leq u,\ 1\leq j\leq p}=\left(  a_{i_{i},j}\right)
_{1\leq i\leq u,\ 1\leq j\leq p}%
\]
\footnote{The double use of the letter \textquotedblleft$i$\textquotedblright%
\ in \textquotedblleft$i_{i}$\textquotedblright\ might appear confusing. The
first \textquotedblleft$i$\textquotedblright\ is part of the notation $i_{k}$
for $k\in\left\{  1,2,\ldots,u\right\}  $; the second \textquotedblleft%
$i$\textquotedblright\ is an element of $\left\{  1,2,\ldots,u\right\}  $.
These two \textquotedblleft$i$\textquotedblright s are unrelated to each
other. I hope the reader can easily tell them apart by the fact that the
\textquotedblleft$i$\textquotedblright\ that is part of the notation $i_{k}$
always appears with a subscript, whereas the second \textquotedblleft%
$i$\textquotedblright\ never does.} (here, we renamed the index $\left(
x,j\right)  $ as $\left(  i,j\right)  $).

Also, $B=\left(  b_{i,j}\right)  _{1\leq i\leq p,\ 1\leq j\leq m}$. Hence, the
definition of $\operatorname*{cols}\nolimits_{j_{1},j_{2},\ldots,j_{v}}B$
yields%
\[
\operatorname*{cols}\nolimits_{j_{1},j_{2},\ldots,j_{v}}B=\left(  b_{i,j_{y}%
}\right)  _{1\leq i\leq p,\ 1\leq y\leq v}=\left(  b_{i,j_{j}}\right)  _{1\leq
i\leq p,\ 1\leq j\leq v}%
\]
(here, we renamed the index $\left(  i,y\right)  $ as $\left(  i,j\right)  $).

We have $\operatorname*{rows}\nolimits_{i_{1},i_{2},\ldots,i_{u}}A=\left(
a_{i_{i},j}\right)  _{1\leq i\leq u,\ 1\leq j\leq p}$ and
$\operatorname*{cols}\nolimits_{j_{1},j_{2},\ldots,j_{v}}B=\left(  b_{i,j_{j}%
}\right)  _{1\leq i\leq p,\ 1\leq j\leq v}$. The definition of the product of
two matrices thus yields%
\begin{align*}
\left(  \operatorname*{rows}\nolimits_{i_{1},i_{2},\ldots,i_{u}}A\right)
\cdot\left(  \operatorname*{cols}\nolimits_{j_{1},j_{2},\ldots,j_{v}}B\right)
&  =\left(  \sum_{k=1}^{p}a_{i_{i},k}b_{k,j_{j}}\right)  _{1\leq i\leq
u,\ 1\leq j\leq v}\\
&  =\left(  \sum_{k=1}^{p}a_{i_{x},k}b_{k,j_{y}}\right)  _{1\leq x\leq
u,\ 1\leq y\leq v}%
\end{align*}
(here, we have renamed the index $\left(  i,j\right)  $ as $\left(
x,y\right)  $). Comparing this with (\ref{pf.lem.adj(AB).minor-of-AB.1}), we
obtain%
\[
\operatorname*{sub}\nolimits_{i_{1},i_{2},\ldots,i_{u}}^{j_{1},j_{2}%
,\ldots,j_{v}}\left(  AB\right)  =\left(  \operatorname*{rows}\nolimits_{i_{1}%
,i_{2},\ldots,i_{u}}A\right)  \cdot\left(  \operatorname*{cols}%
\nolimits_{j_{1},j_{2},\ldots,j_{v}}B\right)  .
\]

\end{verlong}

\noindent This proves Lemma \ref{lem.adj(AB).minor-of-AB}.
\end{proof}

We note in passing that Lemma \ref{lem.adj(AB).minor-of-AB} leads to the
following generalization of Theorem \ref{thm.cauchy-binet}:

\begin{corollary}
\label{cor.adj(AB).cauchy-binet-general}Let $n\in\mathbb{N}$, $m\in\mathbb{N}$
and $p\in\mathbb{N}$. Let $A$ be an $n\times p$-matrix. Let $B$ be a $p\times
m$-matrix. Let $u\in\mathbb{N}$. Let $i_{1},i_{2},\ldots,i_{u}$ be some
elements of $\left\{  1,2,\ldots,n\right\}  $. Let $j_{1},j_{2},\ldots,j_{u}$
be some elements of $\left\{  1,2,\ldots,m\right\}  $. Then,%
\[
\det\left(  \operatorname*{sub}\nolimits_{i_{1},i_{2},\ldots,i_{u}}%
^{j_{1},j_{2},\ldots,j_{u}}\left(  AB\right)  \right)  =\sum_{1\leq
g_{1}<g_{2}<\cdots<g_{u}\leq p}\det\left(  \operatorname*{sub}\nolimits_{i_{1}%
,i_{2},\ldots,i_{u}}^{g_{1},g_{2},\ldots,g_{u}}A\right)  \cdot\det\left(
\operatorname*{sub}\nolimits_{g_{1},g_{2},\ldots,g_{u}}^{j_{1},j_{2}%
,\ldots,j_{u}}B\right)  .
\]
(Here, the summation sign \textquotedblleft$\sum_{1\leq g_{1}<g_{2}%
<\cdots<g_{u}\leq p}$\textquotedblright\ has to be interpreted as
\textquotedblleft$\sum_{\substack{\left(  g_{1},g_{2},\ldots,g_{u}\right)
\in\left\{  1,2,\ldots,p\right\}  ^{u};\\g_{1}<g_{2}<\cdots<g_{u}}%
}$\textquotedblright, in analogy to Remark \ref{rmk.cauchy-binet.sumsign}.)
\end{corollary}

Corollary \ref{cor.adj(AB).cauchy-binet-general} is precisely the formula
\cite[(1.10)]{NoumiYamada}\footnote{We notice that the notation $A_{j_{1}%
,j_{2},\ldots,j_{u}}^{i_{1},i_{2},\ldots,i_{u}}$ in \cite{NoumiYamada} is
equivalent to our notation $\operatorname*{sub}\nolimits_{i_{1},i_{2}%
,\ldots,i_{u}}^{j_{1},j_{2},\ldots,j_{u}}A$.}. We shall use Corollary
\ref{cor.adj(AB).cauchy-binet-general} in a later section, so let us prove it:

\begin{proof}
[Proof of Corollary \ref{cor.adj(AB).cauchy-binet-general}.]Fix any $\left(
g_{1},g_{2},\ldots,g_{u}\right)  \in\left\{  1,2,\ldots,p\right\}  ^{u}$.
Applying Proposition \ref{prop.submatrix.easy} \textbf{(d)} to $p$, $u$ and
$\left(  g_{1},g_{2},\ldots,g_{u}\right)  $ instead of $m$, $v$ and $\left(
j_{1},j_{2},\ldots,j_{v}\right)  $, we obtain%
\[
\operatorname*{sub}\nolimits_{i_{1},i_{2},\ldots,i_{u}}^{g_{1},g_{2}%
,\ldots,g_{u}}A=\operatorname*{rows}\nolimits_{i_{1},i_{2},\ldots,i_{u}%
}\left(  \operatorname*{cols}\nolimits_{g_{1},g_{2},\ldots,g_{u}}A\right)
=\operatorname*{cols}\nolimits_{g_{1},g_{2},\ldots,g_{u}}\left(
\operatorname*{rows}\nolimits_{i_{1},i_{2},\ldots,i_{u}}A\right)  .
\]
Applying Proposition \ref{prop.submatrix.easy} \textbf{(d)} to $p$, $B$, $u$
and $\left(  g_{1},g_{2},\ldots,g_{u}\right)  $ instead of $n$, $A$, $v$ and
$\left(  i_{1},i_{2},\ldots,i_{u}\right)  $, we obtain%
\[
\operatorname*{sub}\nolimits_{g_{1},g_{2},\ldots,g_{u}}^{j_{1},j_{2}%
,\ldots,j_{u}}B=\operatorname*{rows}\nolimits_{g_{1},g_{2},\ldots,g_{u}%
}\left(  \operatorname*{cols}\nolimits_{j_{1},j_{2},\ldots,j_{u}}B\right)
=\operatorname*{cols}\nolimits_{j_{1},j_{2},\ldots,j_{u}}\left(
\operatorname*{rows}\nolimits_{g_{1},g_{2},\ldots,g_{u}}B\right)  .
\]

Now, Lemma \ref{lem.adj(AB).minor-of-AB} (applied to $v=u$) shows that%
\[
\operatorname*{sub}\nolimits_{i_{1},i_{2},\ldots,i_{u}}^{j_{1},j_{2}%
,\ldots,j_{u}}\left(  AB\right)  =\left(  \operatorname*{rows}\nolimits_{i_{1}%
,i_{2},\ldots,i_{u}}A\right)  \cdot\left(  \operatorname*{cols}%
\nolimits_{j_{1},j_{2},\ldots,j_{u}}B\right)  .
\]
Hence,%
\begin{align*}
&  \det\left(  \underbrace{\operatorname*{sub}\nolimits_{i_{1},i_{2}%
,\ldots,i_{u}}^{j_{1},j_{2},\ldots,j_{u}}\left(  AB\right)  }_{=\left(
\operatorname*{rows}\nolimits_{i_{1},i_{2},\ldots,i_{u}}A\right)  \cdot\left(
\operatorname*{cols}\nolimits_{j_{1},j_{2},\ldots,j_{u}}B\right)  }\right) \\
&  =\det\left(  \left(  \operatorname*{rows}\nolimits_{i_{1},i_{2}%
,\ldots,i_{u}}A\right)  \cdot\left(  \operatorname*{cols}\nolimits_{j_{1}%
,j_{2},\ldots,j_{u}}B\right)  \right) \\
&  =\sum_{1\leq g_{1}<g_{2}<\cdots<g_{u}\leq p}\det\left(
\underbrace{\operatorname*{cols}\nolimits_{g_{1},g_{2},\ldots,g_{u}}\left(
\operatorname*{rows}\nolimits_{i_{1},i_{2},\ldots,i_{u}}A\right)
}_{=\operatorname*{sub}\nolimits_{i_{1},i_{2},\ldots,i_{u}}^{g_{1}%
,g_{2},\ldots,g_{u}}A}\right) \\
&  \ \ \ \ \ \ \ \ \ \ \cdot\det\left(  \underbrace{\operatorname*{rows}%
\nolimits_{g_{1},g_{2},\ldots,g_{u}}\left(  \operatorname*{cols}%
\nolimits_{j_{1},j_{2},\ldots,j_{u}}B\right)  }_{=\operatorname*{sub}%
\nolimits_{g_{1},g_{2},\ldots,g_{u}}^{j_{1},j_{2},\ldots,j_{u}}B}\right) \\
&  \ \ \ \ \ \ \ \ \ \ \left(
\begin{array}
[c]{c}%
\text{by Theorem \ref{thm.cauchy-binet} (applied to }u\text{, }p\text{,
}\operatorname*{rows}\nolimits_{i_{1},i_{2},\ldots,i_{u}}A\text{ and}\\
\operatorname*{cols}\nolimits_{j_{1},j_{2},\ldots,j_{u}}B\text{ instead of
}n\text{, }m\text{, }A\text{ and }B\text{)}%
\end{array}
\right) \\
&  =\sum_{1\leq g_{1}<g_{2}<\cdots<g_{u}\leq p}\det\left(  \operatorname*{sub}%
\nolimits_{i_{1},i_{2},\ldots,i_{u}}^{g_{1},g_{2},\ldots,g_{u}}A\right)
\cdot\det\left(  \operatorname*{sub}\nolimits_{g_{1},g_{2},\ldots,g_{u}%
}^{j_{1},j_{2},\ldots,j_{u}}B\right)  .
\end{align*}
This proves Corollary \ref{cor.adj(AB).cauchy-binet-general}.
\end{proof}

\begin{vershort}
\begin{proof}
[Solution to Exercise \ref{exe.adj(AB)}.]Let $\left(  u,v\right)  \in\left\{
1,2,\ldots,n\right\}  ^{2}$. Thus, $1\leq u\leq n$, so that $n\geq1$.

For every $k\in\left\{  1,2,\ldots,n\right\}  $, we have%
\begin{equation}
\operatorname*{cols}\nolimits_{1,2,\ldots,\widehat{k},\ldots,n}\left(
\operatorname*{rows}\nolimits_{1,2,\ldots,\widehat{u},\ldots,n}A\right)
=A_{\sim u,\sim k} \label{sol.adj(AB).short.A}%
\end{equation}
\footnote{\textit{Proof of (\ref{sol.adj(AB).short.A}):} Let $k\in\left\{
1,2,\ldots,n\right\}  $. We can apply Proposition \ref{prop.submatrix.easy}
\textbf{(d)} to $n-1$, $n-1$, $\left(  1,2,\ldots,\widehat{u},\ldots,n\right)
$ and $\left(  1,2,\ldots,\widehat{k},\ldots,n\right)  $ instead of $u$, $v$,
$\left(  i_{1},i_{2},\ldots,i_{u}\right)  $ and $\left(  j_{1},j_{2}%
,\ldots,j_{v}\right)  $. As a result, we obtain%
\[
\operatorname*{sub}\nolimits_{1,2,\ldots,\widehat{u},\ldots,n}^{1,2,\ldots
,\widehat{k},\ldots,n}A=\operatorname*{rows}\nolimits_{1,2,\ldots
,\widehat{u},\ldots,n}\left(  \operatorname*{cols}\nolimits_{1,2,\ldots
,\widehat{k},\ldots,n}A\right)  =\operatorname*{cols}\nolimits_{1,2,\ldots
,\widehat{k},\ldots,n}\left(  \operatorname*{rows}\nolimits_{1,2,\ldots
,\widehat{u},\ldots,n}A\right)  .
\]
But the definition of $A_{\sim u,\sim k}$ yields%
\[
A_{\sim u,\sim k}=\operatorname*{sub}\nolimits_{1,2,\ldots,\widehat{u}%
,\ldots,n}^{1,2,\ldots,\widehat{k},\ldots,n}A=\operatorname*{cols}%
\nolimits_{1,2,\ldots,\widehat{k},\ldots,n}\left(  \operatorname*{rows}%
\nolimits_{1,2,\ldots,\widehat{u},\ldots,n}A\right)  .
\]
This proves (\ref{sol.adj(AB).short.A}).} and%
\begin{equation}
\operatorname*{rows}\nolimits_{1,2,\ldots,\widehat{k},\ldots,n}\left(
\operatorname*{cols}\nolimits_{1,2,\ldots,\widehat{v},\ldots,n}B\right)
=B_{\sim k,\sim v} \label{sol.adj(AB).short.B}%
\end{equation}
\footnote{This holds for similar reasons.}.

We have%
\begin{equation}
\left(  AB\right)  _{\sim u,\sim v}=\left(  \operatorname*{rows}%
\nolimits_{1,2,\ldots,\widehat{u},\ldots,n}A\right)  \cdot\left(
\operatorname*{cols}\nolimits_{1,2,\ldots,\widehat{v},\ldots,n}B\right)
\label{sol.adj(AB).short.1}%
\end{equation}
\footnote{\textit{Proof of (\ref{sol.adj(AB).short.1}):} The definition of
$\left(  AB\right)  _{\sim u,\sim v}$ yields%
\[
\left(  AB\right)  _{\sim u,\sim v}=\operatorname*{sub}\nolimits_{1,2,\ldots
,\widehat{u},\ldots,n}^{1,2,\ldots,\widehat{v},\ldots,n}\left(  AB\right)
=\left(  \operatorname*{rows}\nolimits_{1,2,\ldots,\widehat{u},\ldots
,n}A\right)  \cdot\left(  \operatorname*{cols}\nolimits_{1,2,\ldots
,\widehat{v},\ldots,n}B\right)
\]
(by Lemma \ref{lem.adj(AB).minor-of-AB}, applied to $m=n$, $p=n$, $u=n-1$,
$\left(  i_{1},i_{2},\ldots,i_{u}\right)  =\left(  1,2,\ldots,\widehat{u}%
,\ldots,n\right)  $ and $\left(  j_{1},j_{2},\ldots,j_{v}\right)  =\left(
1,2,\ldots,\widehat{v},\ldots,n\right)  $). This proves
(\ref{sol.adj(AB).short.1}).}. Taking determinants on both sides of this
equation, we obtain%
\begin{align}
&  \det\left(  \left(  AB\right)  _{\sim u,\sim v}\right) \nonumber\\
&  =\det\left(  \left(  \operatorname*{rows}\nolimits_{1,2,\ldots
,\widehat{u},\ldots,n}A\right)  \cdot\left(  \operatorname*{cols}%
\nolimits_{1,2,\ldots,\widehat{v},\ldots,n}B\right)  \right) \nonumber\\
&  =\sum_{k=1}^{n}\det\left(  \underbrace{\operatorname*{cols}%
\nolimits_{1,2,\ldots,\widehat{k},\ldots,n}\left(  \operatorname*{rows}%
\nolimits_{1,2,\ldots,\widehat{u},\ldots,n}A\right)  }_{\substack{=A_{\sim
u,\sim k}\\\text{(by (\ref{sol.adj(AB).short.A}))}}}\right)  \cdot\det\left(
\underbrace{\operatorname*{rows}\nolimits_{1,2,\ldots,\widehat{k},\ldots
,n}\left(  \operatorname*{cols}\nolimits_{1,2,\ldots,\widehat{v},\ldots
,n}B\right)  }_{\substack{=B_{\sim k,\sim v}\\\text{(by
(\ref{sol.adj(AB).short.B}))}}}\right) \nonumber\\
&  \ \ \ \ \ \ \ \ \ \ \left(
\begin{array}
[c]{c}%
\text{by Lemma \ref{lem.adj(AB).cauchy-binet}, applied to }%
\operatorname*{rows}\nolimits_{1,2,\ldots,\widehat{u},\ldots,n}A\\
\text{and }\operatorname*{cols}\nolimits_{1,2,\ldots,\widehat{v},\ldots
,n}B\text{ instead of }A\text{ and }B
\end{array}
\right) \nonumber\\
&  =\sum_{k=1}^{n}\det\left(  A_{\sim u,\sim k}\right)  \cdot\det\left(
B_{\sim k,\sim v}\right)  . \label{sol.adj(AB).short.5}%
\end{align}

Let us now forget that we fixed $\left(  u,v\right)  $. We thus have proven
(\ref{sol.adj(AB).short.5}) for every $\left(  u,v\right)  \in\left\{
1,2,\ldots,n\right\}  ^{2}$.

Now, the definition of $\operatorname*{adj}\left(  AB\right)  $ yields%
\begin{equation}
\operatorname*{adj}\left(  AB\right)  =\left(  \left(  -1\right)  ^{i+j}%
\det\left(  \left(  AB\right)  _{\sim j,\sim i}\right)  \right)  _{1\leq i\leq
n,\ 1\leq j\leq n}. \label{sol.adj(AB).short.L1}%
\end{equation}
But every $\left(  i,j\right)  \in\left\{  1,2,\ldots,n\right\}  ^{2}$
satisfies%
\begin{align*}
&  \left(  -1\right)  ^{i+j}\underbrace{\det\left(  \left(  AB\right)  _{\sim
j,\sim i}\right)  }_{\substack{=\sum_{k=1}^{n}\det\left(  A_{\sim j,\sim
k}\right)  \cdot\det\left(  B_{\sim k,\sim i}\right)  \\\text{(by
(\ref{sol.adj(AB).short.5}), applied to }\left(  u,v\right)  =\left(
j,i\right)  \text{)}}}\\
&  =\left(  -1\right)  ^{i+j}\sum_{k=1}^{n}\det\left(  A_{\sim j,\sim
k}\right)  \cdot\det\left(  B_{\sim k,\sim i}\right) \\
&  =\sum_{k=1}^{n}\underbrace{\left(  -1\right)  ^{i+j}}_{\substack{=\left(
-1\right)  ^{\left(  i+k\right)  +\left(  k+j\right)  }\\\text{(since
}i+j\equiv i+j+2k=\left(  i+k\right)  +\left(  k+j\right)  \operatorname{mod}%
2\text{)}}}\det\left(  A_{\sim j,\sim k}\right)  \cdot\det\left(  B_{\sim
k,\sim i}\right) \\
&  =\sum_{k=1}^{n}\underbrace{\left(  -1\right)  ^{\left(  i+k\right)
+\left(  k+j\right)  }}_{=\left(  -1\right)  ^{i+k}\left(  -1\right)  ^{k+j}%
}\det\left(  A_{\sim j,\sim k}\right)  \cdot\det\left(  B_{\sim k,\sim
i}\right) \\
&  =\sum_{k=1}^{n}\left(  -1\right)  ^{i+k}\det\left(  B_{\sim k,\sim
i}\right)  \cdot\left(  -1\right)  ^{k+j}\det\left(  A_{\sim j,\sim k}\right)
.
\end{align*}
Thus, (\ref{sol.adj(AB).short.L1}) becomes%
\begin{align}
\operatorname*{adj}\left(  AB\right)   &  =\left(  \underbrace{\left(
-1\right)  ^{i+j}\det\left(  \left(  AB\right)  _{\sim j,\sim i}\right)
}_{=\sum_{k=1}^{n}\left(  -1\right)  ^{i+k}\det\left(  B_{\sim k,\sim
i}\right)  \cdot\left(  -1\right)  ^{k+j}\det\left(  A_{\sim j,\sim k}\right)
}\right)  _{1\leq i\leq n,\ 1\leq j\leq n}\nonumber\\
&  =\left(  \sum_{k=1}^{n}\left(  -1\right)  ^{i+k}\det\left(  B_{\sim k,\sim
i}\right)  \cdot\left(  -1\right)  ^{k+j}\det\left(  A_{\sim j,\sim k}\right)
\right)  _{1\leq i\leq n,\ 1\leq j\leq n} \label{sol.adj(AB).short.L2}%
\end{align}

On the other hand, we have $\operatorname*{adj}B=\left(  \left(  -1\right)
^{i+j}\det\left(  B_{\sim j,\sim i}\right)  \right)  _{1\leq i\leq n,\ 1\leq
j\leq n}$ (by the definition of $\operatorname*{adj}B$) and
$\operatorname*{adj}A=\left(  \left(  -1\right)  ^{i+j}\det\left(  A_{\sim
j,\sim i}\right)  \right)  _{1\leq i\leq n,\ 1\leq j\leq n}$ (by the
definition of $\operatorname*{adj}A$). Therefore, the definition of the
product of two matrices shows that%
\[
\operatorname*{adj}B\cdot\operatorname*{adj}A=\left(  \sum_{k=1}^{n}\left(
-1\right)  ^{i+k}\det\left(  B_{\sim k,\sim i}\right)  \cdot\left(  -1\right)
^{k+j}\det\left(  A_{\sim j,\sim k}\right)  \right)  _{1\leq i\leq n,\ 1\leq
j\leq n}.
\]
Compared with (\ref{sol.adj(AB).short.L2}), this yields $\operatorname*{adj}%
\left(  AB\right)  =\operatorname*{adj}B\cdot\operatorname*{adj}A$. This
solves Exercise \ref{exe.adj(AB)}.
\end{proof}
\end{vershort}

\begin{verlong}
\begin{proof}
[Solution to Exercise \ref{exe.adj(AB)}.]Let $\left(  u,v\right)  \in\left\{
1,2,\ldots,n\right\}  ^{2}$. Thus, $u\in\left\{  1,2,\ldots,n\right\}  $ and
$v\in\left\{  1,2,\ldots,n\right\}  $. From $u\in\left\{  1,2,\ldots
,n\right\}  $, we obtain $1\leq u\leq n$, so that $n\geq1$ and therefore
$n-1\in\mathbb{N}$ and $n>0$.

For every $k\in\left\{  1,2,\ldots,n\right\}  $, we have%
\begin{equation}
\operatorname*{cols}\nolimits_{1,2,\ldots,\widehat{k},\ldots,n}\left(
\operatorname*{rows}\nolimits_{1,2,\ldots,\widehat{u},\ldots,n}A\right)
=A_{\sim u,\sim k} \label{sol.adj(AB).A}%
\end{equation}
\footnote{\textit{Proof of (\ref{sol.adj(AB).A}):} Let $k\in\left\{
1,2,\ldots,n\right\}  $. We can apply Proposition \ref{prop.submatrix.easy}
\textbf{(d)} to $n-1$, $n-1$, $\left(  1,2,\ldots,\widehat{u},\ldots,n\right)
$ and $\left(  1,2,\ldots,\widehat{k},\ldots,n\right)  $ instead of $u$, $v$,
$\left(  i_{1},i_{2},\ldots,i_{u}\right)  $ and $\left(  j_{1},j_{2}%
,\ldots,j_{v}\right)  $. As a result, we obtain%
\[
\operatorname*{sub}\nolimits_{1,2,\ldots,\widehat{u},\ldots,n}^{1,2,\ldots
,\widehat{k},\ldots,n}A=\operatorname*{rows}\nolimits_{1,2,\ldots
,\widehat{u},\ldots,n}\left(  \operatorname*{cols}\nolimits_{1,2,\ldots
,\widehat{k},\ldots,n}A\right)  =\operatorname*{cols}\nolimits_{1,2,\ldots
,\widehat{k},\ldots,n}\left(  \operatorname*{rows}\nolimits_{1,2,\ldots
,\widehat{u},\ldots,n}A\right)  .
\]
But the definition of $A_{\sim u,\sim k}$ yields%
\[
A_{\sim u,\sim k}=\operatorname*{sub}\nolimits_{1,2,\ldots,\widehat{u}%
,\ldots,n}^{1,2,\ldots,\widehat{k},\ldots,n}A=\operatorname*{cols}%
\nolimits_{1,2,\ldots,\widehat{k},\ldots,n}\left(  \operatorname*{rows}%
\nolimits_{1,2,\ldots,\widehat{u},\ldots,n}A\right)  .
\]
This proves (\ref{sol.adj(AB).A}).} and%
\begin{equation}
\operatorname*{rows}\nolimits_{1,2,\ldots,\widehat{k},\ldots,n}\left(
\operatorname*{cols}\nolimits_{1,2,\ldots,\widehat{v},\ldots,n}B\right)
=B_{\sim k,\sim v} \label{sol.adj(AB).B}%
\end{equation}
\footnote{\textit{Proof of (\ref{sol.adj(AB).B}):} Let $k\in\left\{
1,2,\ldots,n\right\}  $. We can apply Proposition \ref{prop.submatrix.easy}
\textbf{(d)} to $n-1$, $n-1$, $B$, $\left(  1,2,\ldots,\widehat{k}%
,\ldots,n\right)  $ and $\left(  1,2,\ldots,\widehat{v},\ldots,n\right)  $,
instead of $u$, $v$, $A$, $\left(  i_{1},i_{2},\ldots,i_{u}\right)  $ and
$\left(  j_{1},j_{2},\ldots,j_{v}\right)  $. As a result, we obtain%
\[
\operatorname*{sub}\nolimits_{1,2,\ldots,\widehat{k},\ldots,n}^{1,2,\ldots
,\widehat{v},\ldots,n}B=\operatorname*{rows}\nolimits_{1,2,\ldots
,\widehat{k},\ldots,n}\left(  \operatorname*{cols}\nolimits_{1,2,\ldots
,\widehat{v},\ldots,n}B\right)  =\operatorname*{cols}\nolimits_{1,2,\ldots
,\widehat{v},\ldots,n}\left(  \operatorname*{rows}\nolimits_{1,2,\ldots
,\widehat{k},\ldots,n}B\right)  .
\]
But the definition of $B_{\sim k,\sim v}$ yields%
\[
B_{\sim k,\sim v}=\operatorname*{sub}\nolimits_{1,2,\ldots,\widehat{k}%
,\ldots,n}^{1,2,\ldots,\widehat{v},\ldots,n}B=\operatorname*{rows}%
\nolimits_{1,2,\ldots,\widehat{k},\ldots,n}\left(  \operatorname*{cols}%
\nolimits_{1,2,\ldots,\widehat{v},\ldots,n}B\right)  .
\]
This proves (\ref{sol.adj(AB).B}).}.

We have%
\begin{equation}
\left(  AB\right)  _{\sim u,\sim v}=\left(  \operatorname*{rows}%
\nolimits_{1,2,\ldots,\widehat{u},\ldots,n}A\right)  \cdot\left(
\operatorname*{cols}\nolimits_{1,2,\ldots,\widehat{v},\ldots,n}B\right)
\label{sol.adj(AB).1}%
\end{equation}
\footnote{\textit{Proof of (\ref{sol.adj(AB).1}):} The definition of $\left(
AB\right)  _{\sim u,\sim v}$ yields%
\[
\left(  AB\right)  _{\sim u,\sim v}=\operatorname*{sub}\nolimits_{1,2,\ldots
,\widehat{u},\ldots,n}^{1,2,\ldots,\widehat{v},\ldots,n}\left(  AB\right)
=\left(  \operatorname*{rows}\nolimits_{1,2,\ldots,\widehat{u},\ldots
,n}A\right)  \cdot\left(  \operatorname*{cols}\nolimits_{1,2,\ldots
,\widehat{v},\ldots,n}B\right)
\]
(by Lemma \ref{lem.adj(AB).minor-of-AB}, applied to $m=n$, $p=n$, $u=n-1$,
$\left(  i_{1},i_{2},\ldots,i_{u}\right)  =\left(  1,2,\ldots,\widehat{u}%
,\ldots,n\right)  $ and $\left(  j_{1},j_{2},\ldots,j_{v}\right)  =\left(
1,2,\ldots,\widehat{v},\ldots,n\right)  $). This proves (\ref{sol.adj(AB).1}%
).}. Hence,%
\begin{align}
&  \det\left(  \underbrace{\left(  AB\right)  _{\sim u,\sim v}}_{=\left(
\operatorname*{rows}\nolimits_{1,2,\ldots,\widehat{u},\ldots,n}A\right)
\cdot\left(  \operatorname*{cols}\nolimits_{1,2,\ldots,\widehat{v},\ldots
,n}B\right)  }\right) \nonumber\\
&  =\det\left(  \left(  \operatorname*{rows}\nolimits_{1,2,\ldots
,\widehat{u},\ldots,n}A\right)  \cdot\left(  \operatorname*{cols}%
\nolimits_{1,2,\ldots,\widehat{v},\ldots,n}B\right)  \right) \nonumber\\
&  =\sum_{k=1}^{n}\det\left(  \underbrace{\operatorname*{cols}%
\nolimits_{1,2,\ldots,\widehat{k},\ldots,n}\left(  \operatorname*{rows}%
\nolimits_{1,2,\ldots,\widehat{u},\ldots,n}A\right)  }_{\substack{=A_{\sim
u,\sim k}\\\text{(by (\ref{sol.adj(AB).A}))}}}\right)  \cdot\det\left(
\underbrace{\operatorname*{rows}\nolimits_{1,2,\ldots,\widehat{k},\ldots
,n}\left(  \operatorname*{cols}\nolimits_{1,2,\ldots,\widehat{v},\ldots
,n}B\right)  }_{\substack{=B_{\sim k,\sim v}\\\text{(by (\ref{sol.adj(AB).B}%
))}}}\right) \nonumber\\
&  \ \ \ \ \ \ \ \ \ \ \left(
\begin{array}
[c]{c}%
\text{by Lemma \ref{lem.adj(AB).cauchy-binet}, applied to }%
\operatorname*{rows}\nolimits_{1,2,\ldots,\widehat{u},\ldots,n}A\\
\text{and }\operatorname*{cols}\nolimits_{1,2,\ldots,\widehat{v},\ldots
,n}B\text{ instead of }A\text{ and }B
\end{array}
\right) \nonumber\\
&  =\sum_{k=1}^{n}\det\left(  A_{\sim u,\sim k}\right)  \cdot\det\left(
B_{\sim k,\sim v}\right)  . \label{sol.adj(AB).5}%
\end{align}

Let us now forget that we fixed $\left(  u,v\right)  $. We thus have proven
(\ref{sol.adj(AB).5}) for every $\left(  u,v\right)  \in\left\{
1,2,\ldots,n\right\}  ^{2}$.

Now, the definition of $\operatorname*{adj}\left(  AB\right)  $ yields%
\begin{equation}
\operatorname*{adj}\left(  AB\right)  =\left(  \left(  -1\right)  ^{i+j}%
\det\left(  \left(  AB\right)  _{\sim j,\sim i}\right)  \right)  _{1\leq i\leq
n,\ 1\leq j\leq n}. \label{sol.adj(AB).L1}%
\end{equation}
But every $\left(  i,j\right)  \in\left\{  1,2,\ldots,n\right\}  ^{2}$
satisfies%
\begin{align*}
&  \left(  -1\right)  ^{i+j}\underbrace{\det\left(  \left(  AB\right)  _{\sim
j,\sim i}\right)  }_{\substack{=\sum_{k=1}^{n}\det\left(  A_{\sim j,\sim
k}\right)  \cdot\det\left(  B_{\sim k,\sim i}\right)  \\\text{(by
(\ref{sol.adj(AB).5}), applied to }\left(  u,v\right)  =\left(  j,i\right)
\text{)}}}\\
&  =\left(  -1\right)  ^{i+j}\sum_{k=1}^{n}\det\left(  A_{\sim j,\sim
k}\right)  \cdot\det\left(  B_{\sim k,\sim i}\right) \\
&  =\sum_{k=1}^{n}\underbrace{\left(  -1\right)  ^{i+j}}_{\substack{=\left(
-1\right)  ^{\left(  i+k\right)  +\left(  k+j\right)  }\\\text{(since
}i+j\equiv i+j+2k=\left(  i+k\right)  +\left(  k+j\right)  \operatorname{mod}%
2\text{)}}}\det\left(  A_{\sim j,\sim k}\right)  \cdot\det\left(  B_{\sim
k,\sim i}\right) \\
&  =\sum_{k=1}^{n}\underbrace{\left(  -1\right)  ^{\left(  i+k\right)
+\left(  k+j\right)  }}_{=\left(  -1\right)  ^{i+k}\left(  -1\right)  ^{k+j}%
}\det\left(  A_{\sim j,\sim k}\right)  \cdot\det\left(  B_{\sim k,\sim
i}\right) \\
&  =\sum_{k=1}^{n}\left(  -1\right)  ^{i+k}\underbrace{\left(  -1\right)
^{k+j}\det\left(  A_{\sim j,\sim k}\right)  \cdot\det\left(  B_{\sim k,\sim
i}\right)  }_{=\det\left(  B_{\sim k,\sim i}\right)  \cdot\left(  -1\right)
^{k+j}\det\left(  A_{\sim j,\sim k}\right)  }\\
&  =\sum_{k=1}^{n}\left(  -1\right)  ^{i+k}\det\left(  B_{\sim k,\sim
i}\right)  \cdot\left(  -1\right)  ^{k+j}\det\left(  A_{\sim j,\sim k}\right)
.
\end{align*}
Thus, (\ref{sol.adj(AB).L1}) becomes%
\begin{align}
\operatorname*{adj}\left(  AB\right)   &  =\left(  \underbrace{\left(
-1\right)  ^{i+j}\det\left(  \left(  AB\right)  _{\sim j,\sim i}\right)
}_{=\sum_{k=1}^{n}\left(  -1\right)  ^{i+k}\det\left(  B_{\sim k,\sim
i}\right)  \cdot\left(  -1\right)  ^{k+j}\det\left(  A_{\sim j,\sim k}\right)
}\right)  _{1\leq i\leq n,\ 1\leq j\leq n}\nonumber\\
&  =\left(  \sum_{k=1}^{n}\left(  -1\right)  ^{i+k}\det\left(  B_{\sim k,\sim
i}\right)  \cdot\left(  -1\right)  ^{k+j}\det\left(  A_{\sim j,\sim k}\right)
\right)  _{1\leq i\leq n,\ 1\leq j\leq n} \label{sol.adj(AB).L2}%
\end{align}

On the other hand, we have $\operatorname*{adj}B=\left(  \left(  -1\right)
^{i+j}\det\left(  B_{\sim j,\sim i}\right)  \right)  _{1\leq i\leq n,\ 1\leq
j\leq n}$ (by the definition of $\operatorname*{adj}B$) and
$\operatorname*{adj}A=\left(  \left(  -1\right)  ^{i+j}\det\left(  A_{\sim
j,\sim i}\right)  \right)  _{1\leq i\leq n,\ 1\leq j\leq n}$ (by the
definition of $\operatorname*{adj}A$). Therefore, the definition of the
product of two matrices shows that%
\[
\operatorname*{adj}B\cdot\operatorname*{adj}A=\left(  \sum_{k=1}^{n}\left(
-1\right)  ^{i+k}\det\left(  B_{\sim k,\sim i}\right)  \cdot\left(  -1\right)
^{k+j}\det\left(  A_{\sim j,\sim k}\right)  \right)  _{1\leq i\leq n,\ 1\leq
j\leq n}.
\]
Compared with (\ref{sol.adj(AB).L2}), this yields $\operatorname*{adj}\left(
AB\right)  =\operatorname*{adj}B\cdot\operatorname*{adj}A$. This solves
Exercise \ref{exe.adj(AB)}.
\end{proof}
\end{verlong}

\subsection{Solution to Exercise \ref{exe.vander-hook}}

Throughout this section, we shall use the notation $V\left(  y_{1}%
,y_{2},\ldots,y_{n}\right)  $ defined in Exercise \ref{exe.vander-hook}. Let
us now make some preparations for solving Exercise \ref{exe.vander-hook}.

\subsubsection{Lemmas}

\begin{definition}
If $i$ and $j$ are two objects, then $\delta_{i,j}$ is defined to be the
element $%
\begin{cases}
1, & \text{if }i=j;\\
0, & \text{if }i\neq j
\end{cases}
$ of $\mathbb{K}$.
\end{definition}

Let us now prove four mostly trivial lemmas:

\begin{lemma}
\label{lem.sol.vander-hook.gauss}Let $n\in\mathbb{N}$. Then,%
\begin{equation}
\sum_{r=0}^{n-1}r=\dbinom{n}{2}. \label{sol.vander-hook.gauss}%
\end{equation}

\end{lemma}

\begin{vershort}
\begin{proof}
[Proof of Lemma \ref{lem.sol.vander-hook.gauss}.]If $n=0$, then Lemma
\ref{lem.sol.vander-hook.gauss} holds for obvious reasons. Thus, we WLOG
assume that $n\neq0$. Hence, $n\geq1$, so that $n-1\in\mathbb{N}$. Now,%
\begin{align*}
\sum_{r=0}^{n-1}r  &  =\sum_{i=0}^{n-1}i\ \ \ \ \ \ \ \ \ \ \left(
\text{here, we have renamed the summation index }r\text{ as }i\right) \\
&  =\dfrac{\left(  n-1\right)  \left(  \left(  n-1\right)  +1\right)  }%
{2}\ \ \ \ \ \ \ \ \ \ \left(  \text{by (\ref{eq.sum.littlegauss1}) (applied
to }n-1\text{ instead of }n\text{)}\right) \\
&  =\dfrac{n\left(  n-1\right)  }{2}=\dbinom{n}{2}.
\end{align*}
This proves Lemma \ref{lem.sol.vander-hook.gauss}.
\end{proof}
\end{vershort}

\begin{verlong}
\begin{proof}
[Proof of Lemma \ref{lem.sol.vander-hook.gauss}.]If $n=0$, then%
\begin{align*}
\sum_{r=0}^{n-1}r  &  =\sum_{r=0}^{0-1}r=\left(  \text{empty sum}\right)
=0=\dbinom{0}{2}\ \ \ \ \ \ \ \ \ \ \left(  \text{since }\dbinom{0}{2}%
=\dfrac{0\left(  0-1\right)  }{2}=0\right) \\
&  =\dbinom{n}{2}\ \ \ \ \ \ \ \ \ \ \left(  \text{since }0=n\right)  .
\end{align*}
Hence, if $n=0$, then Lemma \ref{lem.sol.vander-hook.gauss} holds. Thus, for
the rest of the proof of Lemma \ref{lem.sol.vander-hook.gauss}, we can WLOG
assume that we don't have $n=0$. Assume this.

We have $n\neq0$ (since we don't have $n=0$). Thus, $n\geq1$ (since
$n\in\mathbb{N}$). Therefore, $n-1\in\mathbb{N}$. Now,%
\begin{align*}
\sum_{r=0}^{n-1}r  &  =\sum_{i=0}^{n-1}i\ \ \ \ \ \ \ \ \ \ \left(
\text{here, we have renamed the summation index }r\text{ as }i\right) \\
&  =\dfrac{\left(  n-1\right)  \left(  \left(  n-1\right)  +1\right)  }%
{2}\ \ \ \ \ \ \ \ \ \ \left(  \text{by (\ref{eq.sum.littlegauss1}) (applied
to }n-1\text{ instead of }n\text{)}\right) \\
&  =\dfrac{\left(  n-1\right)  n}{2}=\dfrac{n\left(  n-1\right)  }{2}%
=\dbinom{n}{2}.
\end{align*}
This proves Lemma \ref{lem.sol.vander-hook.gauss}.
\end{proof}
\end{verlong}

\begin{lemma}
\label{lem.sol.vander-hook.deltas}Let $n\in\mathbb{N}$. Let $x_{1}%
,x_{2},\ldots,x_{n}$ be $n$ elements of $\mathbb{K}$. Let $t\in\mathbb{K}$.
Let $k\in\left\{  1,2,\ldots,n\right\}  $. For each $j\in\left\{
1,2,\ldots,n\right\}  $, set $y_{j}=x_{j}+\delta_{j,k}t$. Then,%
\[
\left(  x_{1},x_{2},\ldots,x_{k-1},x_{k}+t,x_{k+1},x_{k+2},\ldots
,x_{n}\right)  =\left(  y_{1},y_{2},\ldots,y_{n}\right)  .
\]

\end{lemma}

\begin{vershort}
\begin{proof}
[Proof of Lemma \ref{lem.sol.vander-hook.deltas}.]For each $j\in\left\{
1,2,\ldots,n\right\}  $, we have%
\begin{align*}
&  \left(  \text{the }j\text{-th entry of the list }\left(  y_{1},y_{2}%
,\ldots,y_{n}\right)  \right) \\
&  =y_{j}=x_{j}+\underbrace{\delta_{j,k}}_{=%
\begin{cases}
1, & \text{if }j=k;\\
0, & \text{if }j\neq k
\end{cases}
}t=x_{j}+%
\begin{cases}
1, & \text{if }j=k;\\
0, & \text{if }j\neq k
\end{cases}
t\\
&  =%
\begin{cases}
x_{j}+1t, & \text{if }j=k;\\
x_{j}+0t, & \text{if }j\neq k
\end{cases}
=%
\begin{cases}
x_{j}+t, & \text{if }j=k;\\
x_{j}, & \text{if }j\neq k
\end{cases}
\\
&  =%
\begin{cases}
x_{k}+t, & \text{if }j=k;\\
x_{j}, & \text{if }j\neq k
\end{cases}
\ \ \ \ \ \ \ \ \ \ \left(  \text{since }x_{j}=x_{k}\text{ in the case when
}j=k\right) \\
&  =\left(  \text{the }j\text{-th entry of the list }\left(  x_{1}%
,x_{2},\ldots,x_{k-1},x_{k}+t,x_{k+1},x_{k+2},\ldots,x_{n}\right)  \right)  .
\end{align*}
Hence, $\left(  y_{1},y_{2},\ldots,y_{n}\right)  =\left(  x_{1},x_{2}%
,\ldots,x_{k-1},x_{k}+t,x_{k+1},x_{k+2},\ldots,x_{n}\right)  $. This proves
Lemma \ref{lem.sol.vander-hook.deltas}.
\end{proof}
\end{vershort}

\begin{verlong}
\begin{proof}
[Proof of Lemma \ref{lem.sol.vander-hook.deltas}.]If $\alpha$ and $\beta$ are
two finite lists of elements of $\mathbb{K}$, then we define a new finite list
$\alpha\ast\beta$ of elements of $\mathbb{K}$ by setting%
\[
\alpha\ast\beta=\left(  \alpha_{1},\alpha_{2},\ldots,\alpha_{a},\beta
_{1},\beta_{2},\ldots,\beta_{b}\right)  ,
\]
where $\alpha$ is written in the form $\left(  \alpha_{1},\alpha_{2}%
,\ldots,\alpha_{a}\right)  $, and where $\beta$ is written in the form
$\left(  \beta_{1},\beta_{2},\ldots,\beta_{b}\right)  $. This new list
$\alpha\ast\beta$ is called the \textit{concatenation} of the lists $\alpha$
and $\beta$. If $\alpha$, $\beta$ and $\gamma$ are three finite lists of
elements of $\mathbb{K}$, then the two lists $\left(  \alpha\ast\beta\right)
\ast\gamma$ and $\alpha\ast\left(  \beta\ast\gamma\right)  $ are identical. In
other words, concatenation of finite lists is an associative operation. Thus,
we are able to write $\alpha\ast\beta\ast\gamma$ for any of the two identical
lists $\left(  \alpha\ast\beta\right)  \ast\gamma$ and $\alpha\ast\left(
\beta\ast\gamma\right)  $ whenever $\alpha$, $\beta$ and $\gamma$ are three
finite lists of elements of $\mathbb{K}$.

Now, $x_{j}=y_{j}$ for each $j\in\left\{  1,2,\ldots,k-1\right\}
$\ \ \ \ \footnote{\textit{Proof.} Let $j\in\left\{  1,2,\ldots,k-1\right\}
$. Then, $j\leq k-1<k$, hence $j\neq k$. Now, the definition of $y_{j}$ yields
$y_{j}=x_{j}+\underbrace{\delta_{j,k}}_{\substack{=0\\\text{(since }j\neq
k\text{)}}}t=x_{j}+\underbrace{0t}_{=0}=x_{j}$. In other words, $x_{j}=y_{j}$.
Qed.}. In other words, $\left(  x_{1},x_{2},\ldots,x_{k-1}\right)  =\left(
y_{1},y_{2},\ldots,y_{k-1}\right)  $.

The definition of $y_{k}$ yields $y_{k}=x_{k}+\underbrace{\delta_{k,k}%
}_{\substack{=1\\\text{(since }k=k\text{)}}}t=x_{k}+\underbrace{1t}_{=t}%
=x_{k}+t$. In other words, $x_{k}+t=y_{k}$.

Also, $x_{j}=y_{j}$ for each $j\in\left\{  k+1,k+2,\ldots,n\right\}
$\ \ \ \ \footnote{\textit{Proof.} Let $j\in\left\{  k+1,k+2,\ldots,n\right\}
$. Then, $j\geq k+1>k$, hence $j\neq k$. Now, the definition of $y_{j}$ yields
$y_{j}=x_{j}+\underbrace{\delta_{j,k}}_{\substack{=0\\\text{(since }j\neq
k\text{)}}}t=x_{j}+\underbrace{0t}_{=0}=x_{j}$. In other words, $x_{j}=y_{j}$.
Qed.}. In other words, $\left(  x_{k+1},x_{k+2},\ldots,x_{n}\right)  =\left(
y_{k+1},y_{k+2},\ldots,y_{n}\right)  $.

Now,%
\begin{align*}
&  \left(  x_{1},x_{2},\ldots,x_{k-1},x_{k}+t,x_{k+1},x_{k+2},\ldots
,x_{n}\right) \\
&  =\underbrace{\left(  x_{1},x_{2},\ldots,x_{k-1},x_{k}+t\right)  }_{=\left(
x_{1},x_{2},\ldots,x_{k-1}\right)  \ast\left(  x_{k}+t\right)  }\ast\left(
x_{k+1},x_{k+2},\ldots,x_{n}\right) \\
&  =\underbrace{\left(  x_{1},x_{2},\ldots,x_{k-1}\right)  }_{=\left(
y_{1},y_{2},\ldots,y_{k-1}\right)  }\ast\left(  \underbrace{x_{k}+t}_{=y_{k}%
}\right)  \ast\underbrace{\left(  x_{k+1},x_{k+2},\ldots,x_{n}\right)
}_{=\left(  y_{k+1},y_{k+2},\ldots,y_{n}\right)  }\\
&  =\underbrace{\left(  y_{1},y_{2},\ldots,y_{k-1}\right)  \ast\left(
y_{k}\right)  }_{=\left(  y_{1},y_{2},\ldots,y_{k-1},y_{k}\right)  =\left(
y_{1},y_{2},\ldots,y_{k}\right)  }\ast\left(  y_{k+1},y_{k+2},\ldots
,y_{n}\right) \\
&  =\left(  y_{1},y_{2},\ldots,y_{k}\right)  \ast\left(  y_{k+1}%
,y_{k+2},\ldots,y_{n}\right)  =\left(  y_{1},y_{2},\ldots,y_{k},y_{k+1}%
,y_{k+2},\ldots,y_{n}\right) \\
&  =\left(  y_{1},y_{2},\ldots,y_{n}\right)  .
\end{align*}
This proves Lemma \ref{lem.sol.vander-hook.deltas}.
\end{proof}
\end{verlong}

\begin{lemma}
\label{lem.sol.vander-hook.V=}Let $n\in\mathbb{N}$. Let $y_{1},y_{2}%
,\ldots,y_{n}$ be $n$ elements of $\mathbb{K}$. Then,%
\[
V\left(  y_{1},y_{2},\ldots,y_{n}\right)  =\det\left(  \left(  y_{i}%
^{n-j}\right)  _{1\leq i\leq n,\ 1\leq j\leq n}\right)  .
\]

\end{lemma}

\begin{proof}
[Proof of Lemma \ref{lem.sol.vander-hook.V=}.]Theorem \ref{thm.vander-det}
\textbf{(a)} (applied to $x_{i}=y_{i}$) yields%
\[
\det\left(  \left(  y_{i}^{n-j}\right)  _{1\leq i\leq n,\ 1\leq j\leq
n}\right)  =\prod_{1\leq i<j\leq n}\left(  y_{i}-y_{j}\right)  .
\]
Comparing this with%
\[
V\left(  y_{1},y_{2},\ldots,y_{n}\right)  =\prod_{1\leq i<j\leq n}\left(
y_{i}-y_{j}\right)  \ \ \ \ \ \ \ \ \ \ \left(  \text{by the definition of
}V\left(  y_{1},y_{2},\ldots,y_{n}\right)  \right)  ,
\]
we obtain $V\left(  y_{1},y_{2},\ldots,y_{n}\right)  =\det\left(  \left(
y_{i}^{n-j}\right)  _{1\leq i\leq n,\ 1\leq j\leq n}\right)  $. Thus, Lemma
\ref{lem.sol.vander-hook.V=} is proven.
\end{proof}

\begin{lemma}
\label{lem.sol.vander-hook.AvsB}Let $n\in\mathbb{N}$ and $m\in\mathbb{N}$. Let
$k\in\left\{  1,2,\ldots,n\right\}  $. Let $A=\left(  a_{i,j}\right)  _{1\leq
i\leq n,\ 1\leq j\leq m}$ and $B=\left(  b_{i,j}\right)  _{1\leq i\leq
n,\ 1\leq j\leq m}$ be two $n\times m$-matrices. Assume that every
$i\in\left\{  1,2,\ldots,n\right\}  $ and $j\in\left\{  1,2,\ldots,m\right\}
$ satisfying $i\neq k$ satisfy%
\begin{equation}
a_{i,j}=b_{i,j}. \label{eq.lem.sol.vander-hook.AvsB.ass}%
\end{equation}
Then, $B_{\sim k,\sim q}=A_{\sim k,\sim q}$ for every $q\in\left\{
1,2,\ldots,m\right\}  $.
\end{lemma}

\begin{vershort}
\begin{proof}
[Proof of Lemma \ref{lem.sol.vander-hook.AvsB}.]Let $q\in\left\{
1,2,\ldots,m\right\}  $. The condition (\ref{eq.lem.sol.vander-hook.AvsB.ass})
shows that the matrices $A$ and $B$ agree in all their entries except for
those in the respective $k$-th rows of the matrices. Therefore, the matrices
$A_{\sim k,\sim q}$ and $B_{\sim k,\sim q}$ must agree completely (since they
are obtained from the matrices $A$ and $B$ by removing the $k$-th rows and the
$q$-th columns). In other words, $A_{\sim k,\sim q}=B_{\sim k,\sim q}$. This
proves Lemma \ref{lem.sol.vander-hook.AvsB}.
\end{proof}
\end{vershort}

\begin{verlong}
\begin{proof}
[Proof of Lemma \ref{lem.sol.vander-hook.AvsB}.]Let $q\in\left\{
1,2,\ldots,m\right\}  $.

Denote the list $\left(  1,2,\ldots,\widehat{k},\ldots,n\right)  $ by $\left(
u_{1},u_{2},\ldots,u_{n-1}\right)  $. (This is possible, since the list
$\left(  1,2,\ldots,\widehat{k},\ldots,n\right)  $ has size $n-1$.) Then,
$\left(  u_{1},u_{2},\ldots,u_{n-1}\right)  =\left(  1,2,\ldots,\widehat{k}%
,\ldots,n\right)  $. Now, every $x\in\left\{  1,2,\ldots,n-1\right\}  $
satisfies%
\begin{align}
u_{x}  &  \in\left\{  u_{1},u_{2},\ldots,u_{n-1}\right\}  =\left\{
1,2,\ldots,\widehat{k},\ldots,n\right\} \nonumber\\
&  \ \ \ \ \ \ \ \ \ \ \left(  \text{since }\left(  u_{1},u_{2},\ldots
,u_{n-1}\right)  =\left(  1,2,\ldots,\widehat{k},\ldots,n\right)  \right)
\nonumber\\
&  =\left\{  1,2,\ldots,n\right\}  \setminus\left\{  k\right\}  .
\label{pf.lem.sol.vander-hook.AvsB.uxin}%
\end{align}

Denote the list $\left(  1,2,\ldots,\widehat{q},\ldots,m\right)  $ by $\left(
v_{1},v_{2},\ldots,v_{m-1}\right)  $. (This is possible, since the list
$\left(  1,2,\ldots,\widehat{q},\ldots,m\right)  $ has size $m-1$.) Then,
$\left(  v_{1},v_{2},\ldots,v_{m-1}\right)  =\left(  1,2,\ldots,\widehat{q}%
,\ldots,m\right)  $. Now, every $y\in\left\{  1,2,\ldots,m-1\right\}  $
satisfies%
\begin{align}
v_{y}  &  \in\left\{  v_{1},v_{2},\ldots,v_{m-1}\right\}  =\left\{
1,2,\ldots,\widehat{q},\ldots,m\right\} \nonumber\\
&  \ \ \ \ \ \ \ \ \ \ \left(  \text{since }\left(  v_{1},v_{2},\ldots
,v_{m-1}\right)  =\left(  1,2,\ldots,\widehat{q},\ldots,m\right)  \right)
\nonumber\\
&  =\left\{  1,2,\ldots,m\right\}  \setminus\left\{  q\right\}  .
\label{pf.lem.sol.vander-hook.AvsB.vyin}%
\end{align}

Every $\left(  x,y\right)  \in\left\{  1,2,\ldots,n-1\right\}  \times\left\{
1,2,\ldots,m-1\right\}  $ satisfies%
\begin{equation}
a_{u_{x},v_{y}}=b_{u_{x},v_{y}} \label{pf.lem.sol.vander-hook.AvsB.a=b}%
\end{equation}
\footnote{\textit{Proof of (\ref{pf.lem.sol.vander-hook.AvsB.a=b}):} Let
$\left(  x,y\right)  \in\left\{  1,2,\ldots,n-1\right\}  \times\left\{
1,2,\ldots,m-1\right\}  $. Thus, $x\in\left\{  1,2,\ldots,n-1\right\}  $ and
$y\in\left\{  1,2,\ldots,m-1\right\}  $.
\par
Now, (\ref{pf.lem.sol.vander-hook.AvsB.uxin}) shows that $u_{x}\in\left\{
1,2,\ldots,n\right\}  \setminus\left\{  k\right\}  $. In other words,
$u_{x}\in\left\{  1,2,\ldots,n\right\}  $ and $u_{x}\neq k$. Also,
(\ref{pf.lem.sol.vander-hook.AvsB.vyin}) yields $v_{y}\in\left\{
1,2,\ldots,m\right\}  \setminus\left\{  q\right\}  \subseteq\left\{
1,2,\ldots,m\right\}  $. Hence, (\ref{eq.lem.sol.vander-hook.AvsB.ass})
(applied to $i=u_{x}$ and $j=v_{y}$) yields $a_{u_{x},v_{y}}=b_{u_{x},v_{y}}$.
This proves (\ref{pf.lem.sol.vander-hook.AvsB.a=b}).}.

The definition of $A_{\sim k,\sim q}$ yields
\begin{align}
A_{\sim k,\sim q}  &  =\operatorname*{sub}\nolimits_{1,2,\ldots,\widehat{k}%
,\ldots,n}^{1,2,\ldots,\widehat{q},\ldots,m}A=\operatorname*{sub}%
\nolimits_{u_{1},u_{2},\ldots,u_{n-1}}^{1,2,\ldots,\widehat{q},\ldots
,m}A\nonumber\\
&  \ \ \ \ \ \ \ \ \ \ \left(  \text{since }\left(  1,2,\ldots,\widehat{k}%
,\ldots,n\right)  =\left(  u_{1},u_{2},\ldots,u_{n-1}\right)  \right)
\nonumber\\
&  =\operatorname*{sub}\nolimits_{u_{1},u_{2},\ldots,u_{n-1}}^{v_{1}%
,v_{2},\ldots,v_{m-1}}A\nonumber\\
&  \ \ \ \ \ \ \ \ \ \ \left(  \text{since }\left(  1,2,\ldots,\widehat{q}%
,\ldots,m\right)  =\left(  v_{1},v_{2},\ldots,v_{m-1}\right)  \right)
\nonumber\\
&  =\left(  a_{u_{x},v_{y}}\right)  _{1\leq x\leq n-1,\ 1\leq y\leq m-1}
\label{pf.lem.sol.vander-hook.AvsB.1A}%
\end{align}
(by the definition of $\operatorname*{sub}\nolimits_{u_{1},u_{2}%
,\ldots,u_{n-1}}^{v_{1},v_{2},\ldots,v_{m-1}}A$, since $A=\left(
a_{i,j}\right)  _{1\leq i\leq n,\ 1\leq j\leq m}$). The same argument (but
applied to $B$ and $b_{i,j}$ instead of $A$ and $a_{i,j}$) yields%
\begin{equation}
B_{\sim k,\sim q}=\left(  b_{u_{x},v_{y}}\right)  _{1\leq x\leq n-1,\ 1\leq
y\leq m-1}. \label{pf.lem.sol.vander-hook.AvsB.1B}%
\end{equation}
Comparing this with%
\begin{align*}
A_{\sim k,\sim q}  &  =\left(  \underbrace{a_{u_{x},v_{y}}}%
_{\substack{=b_{u_{x},v_{y}}\\\text{(by (\ref{pf.lem.sol.vander-hook.AvsB.a=b}%
))}}}\right)  _{1\leq x\leq n-1,\ 1\leq y\leq m-1}\ \ \ \ \ \ \ \ \ \ \left(
\text{by (\ref{pf.lem.sol.vander-hook.AvsB.1A})}\right) \\
&  =\left(  b_{u_{x},v_{y}}\right)  _{1\leq x\leq n-1,\ 1\leq y\leq m-1},
\end{align*}
we obtain $B_{\sim k,\sim q}=A_{\sim k,\sim q}$. This proves Lemma
\ref{lem.sol.vander-hook.AvsB}.
\end{proof}
\end{verlong}

Next, let us combine Theorem \ref{thm.laplace.gen} \textbf{(b)} and
Proposition \ref{prop.laplace.0} \textbf{(b)} into a convenient package:

\begin{lemma}
\label{lem.sol.vander-hook.lap0}Let $n\in\mathbb{N}$. Let $A=\left(
a_{i,j}\right)  _{1\leq i\leq n,\ 1\leq j\leq n}$ be an $n\times n$-matrix.
Let $q\in\left\{  1,2,\ldots,n\right\}  $ and $r\in\left\{  1,2,\ldots
,n\right\}  $. Then,%
\[
\sum_{p=1}^{n}\left(  -1\right)  ^{p+q}a_{p,r}\det\left(  A_{\sim p,\sim
q}\right)  =\delta_{q,r}\det A.
\]

\end{lemma}

\begin{proof}
[Proof of Lemma \ref{lem.sol.vander-hook.lap0}.]We are in one of the following
two cases:

\textit{Case 1:} We have $q=r$.

\textit{Case 2:} We have $q\neq r$.

Let us first consider Case 1. In this case, we have $q=r$. Hence,
$\delta_{q,r}=1$, so that
\begin{align*}
\underbrace{\delta_{q,r}}_{=1}\det A  &  =\det A=\sum_{p=1}^{n}\left(
-1\right)  ^{p+q}\underbrace{a_{p,q}}_{\substack{=a_{p,r}\\\text{(since
}q=r\text{)}}}\det\left(  A_{\sim p,\sim q}\right)
\ \ \ \ \ \ \ \ \ \ \left(  \text{by Theorem \ref{thm.laplace.gen}
\textbf{(b)}}\right) \\
&  =\sum_{p=1}^{n}\left(  -1\right)  ^{p+q}a_{p,r}\det\left(  A_{\sim p,\sim
q}\right)  .
\end{align*}
Hence, Lemma \ref{lem.sol.vander-hook.lap0} is proven in Case 1.

Let us now consider Case 2. In this case, we have $q\neq r$. Hence,
$\delta_{q,r}=0$, so that
\begin{align*}
\underbrace{\delta_{q,r}}_{=0}\det A  &  =0\det A=0\\
&  =\sum_{p=1}^{n}\left(  -1\right)  ^{p+q}a_{p,r}\det\left(  A_{\sim p,\sim
q}\right)  \ \ \ \ \ \ \ \ \ \ \left(  \text{by Proposition
\ref{prop.laplace.0} \textbf{(b)}}\right)  .
\end{align*}
Hence, Lemma \ref{lem.sol.vander-hook.lap0} is proven in Case 2.

We have now proven Lemma \ref{lem.sol.vander-hook.lap0} in each of the two
Cases 1 and 2. Since these two Cases cover all possibilities, this shows that
Lemma \ref{lem.sol.vander-hook.lap0} always holds.
\end{proof}

Next, we show three crucial lemmas:

\begin{lemma}
\label{lem.sol.vander-hook.V=x}Let $n\in\mathbb{N}$. Let $x_{1},x_{2}%
,\ldots,x_{n}$ be $n$ elements of $\mathbb{K}$. Let $A$ be the $n\times
n$-matrix $\left(  x_{i}^{n-j}\right)  _{1\leq i\leq n,\ 1\leq j\leq n}$.
Then,%
\begin{equation}
V\left(  x_{1},x_{2},\ldots,x_{n}\right)  =\det A.
\label{eq.lem.sol.vander-hook.V=x.eq}%
\end{equation}

\end{lemma}

\begin{proof}
[Proof of Lemma \ref{lem.sol.vander-hook.V=x}.]The definition of $A$ yields
$A=\left(  x_{i}^{n-j}\right)  _{1\leq i\leq n,\ 1\leq j\leq n}$. Lemma
\ref{lem.sol.vander-hook.V=} (applied to $y_{j}=x_{j}$) yields%
\[
V\left(  x_{1},x_{2},\ldots,x_{n}\right)  =\det\left(  \underbrace{\left(
x_{i}^{n-j}\right)  _{1\leq i\leq n,\ 1\leq j\leq n}}_{=A}\right)  =\det A.
\]
This proves Lemma \ref{lem.sol.vander-hook.V=x}.
\end{proof}

\begin{lemma}
\label{lem.sol.vander-hook.lap1}Let $n\in\mathbb{N}$. Let $x_{1},x_{2}%
,\ldots,x_{n}$ be $n$ elements of $\mathbb{K}$. Let $t\in\mathbb{K}$. Let $A$
be the $n\times n$-matrix $\left(  x_{i}^{n-j}\right)  _{1\leq i\leq n,\ 1\leq
j\leq n}$. Let $k\in\left\{  1,2,\ldots,n\right\}  $.

\textbf{(a)} We have%
\[
V\left(  x_{1},x_{2},\ldots,x_{n}\right)  =\sum_{q=1}^{n}\left(  -1\right)
^{k+q}x_{k}^{n-q}\det\left(  A_{\sim k,\sim q}\right)  .
\]

\textbf{(b)} We have%
\begin{align*}
&  V\left(  x_{1},x_{2},\ldots,x_{k-1},x_{k}+t,x_{k+1},x_{k+2},\ldots
,x_{n}\right) \\
&  =\sum_{q=1}^{n}\left(  -1\right)  ^{k+q}\left(  x_{k}+t\right)  ^{n-q}%
\det\left(  A_{\sim k,\sim q}\right)  .
\end{align*}

\textbf{(c)} We have%
\begin{align*}
&  V\left(  x_{1},x_{2},\ldots,x_{k-1},x_{k}+t,x_{k+1},x_{k+2},\ldots
,x_{n}\right)  -V\left(  x_{1},x_{2},\ldots,x_{n}\right) \\
&  =\sum_{q=1}^{n}\ \ \sum_{\ell=1}^{n-q}\dbinom{n-q}{\ell}t^{\ell}\left(
-1\right)  ^{k+q}x_{k}^{n-q-\ell}\det\left(  A_{\sim k,\sim q}\right)  .
\end{align*}

\end{lemma}

\begin{proof}
[Proof of Lemma \ref{lem.sol.vander-hook.lap1}.]The definition of $A$ yields
$A=\left(  x_{i}^{n-j}\right)  _{1\leq i\leq n,\ 1\leq j\leq n}$.

\textbf{(a)} Lemma \ref{lem.sol.vander-hook.V=x} yields%
\[
V\left(  x_{1},x_{2},\ldots,x_{n}\right)  =\det A=\sum_{q=1}^{n}\left(
-1\right)  ^{k+q}x_{k}^{n-q}\det\left(  A_{\sim k,\sim q}\right)
\]
(by Theorem \ref{thm.laplace.gen} \textbf{(a)}, applied to $x_{i}^{n-j}$ and
$k$ instead of $a_{i,j}$ and $p$). This proves Lemma
\ref{lem.sol.vander-hook.lap1} \textbf{(a)}.

\textbf{(b)} For each $j\in\left\{  1,2,\ldots,n\right\}  $, define an element
$y_{j}\in\mathbb{K}$ by $y_{j}=x_{j}+\delta_{j,k}t$. Define an $n\times
n$-matrix $B$ by $B=\left(  y_{i}^{n-j}\right)  _{1\leq i\leq n,\ 1\leq j\leq
n}$.

Every $i\in\left\{  1,2,\ldots,n\right\}  $ and $j\in\left\{  1,2,\ldots
,n\right\}  $ satisfying $i\neq k$ satisfy $x_{i}^{n-j}=y_{i}^{n-j}%
$\ \ \ \ \footnote{\textit{Proof.} Let $i\in\left\{  1,2,\ldots,n\right\}  $
and $j\in\left\{  1,2,\ldots,n\right\}  $ be such that $i\neq k$. Then, the
definition of $y_{i}$ yields $y_{i}=x_{i}+\underbrace{\delta_{i,k}%
}_{\substack{=0\\\text{(since }i\neq k\text{)}}}t=x_{i}+\underbrace{0t}%
_{=0}=x_{i}$. Hence, $x_{i}=y_{i}$, and thus $x_{i}^{n-j}=y_{i}^{n-j}$. Qed.}.
Hence, Lemma \ref{lem.sol.vander-hook.AvsB} (applied to $n$, $x_{i}^{n-j}$ and
$y_{i}^{n-j}$ instead of $m$, $a_{i,j}$ and $b_{i,j}$) yields that%
\begin{equation}
B_{\sim k,\sim q}=A_{\sim k,\sim q}\ \ \ \ \ \ \ \ \ \ \text{for every }%
q\in\left\{  1,2,\ldots,n\right\}  . \label{pf.lem.sol.vander-hook.lap1.B.1}%
\end{equation}
Moreover, the definition of $y_{k}$ satisfies $y_{k}=x_{k}+\underbrace{\delta
_{k,k}}_{\substack{=1\\\text{(since }k=k\text{)}}}t=x_{k}+\underbrace{1t}%
_{=t}=x_{k}+t$.

Lemma \ref{lem.sol.vander-hook.deltas} yields%
\[
\left(  x_{1},x_{2},\ldots,x_{k-1},x_{k}+t,x_{k+1},x_{k+2},\ldots
,x_{n}\right)  =\left(  y_{1},y_{2},\ldots,y_{n}\right)  .
\]
Hence,%
\begin{align*}
&  V\underbrace{\left(  x_{1},x_{2},\ldots,x_{k-1},x_{k}+t,x_{k+1}%
,x_{k+2},\ldots,x_{n}\right)  }_{=\left(  y_{1},y_{2},\ldots,y_{n}\right)  }\\
&  =V\left(  y_{1},y_{2},\ldots,y_{n}\right)  =\det\left(  \underbrace{\left(
y_{i}^{n-j}\right)  _{1\leq i\leq n,\ 1\leq j\leq n}}_{=B}\right)
\ \ \ \ \ \ \ \ \ \ \left(  \text{by Lemma \ref{lem.sol.vander-hook.V=}%
}\right) \\
&  =\det B=\sum_{q=1}^{n}\left(  -1\right)  ^{k+q}\underbrace{y_{k}^{n-q}%
}_{\substack{=\left(  x_{k}+t\right)  ^{n-q}\\\text{(since }y_{k}%
=x_{k}+t\text{)}}}\det\left(  \underbrace{B_{\sim k,\sim q}}%
_{\substack{=A_{\sim k,\sim q}\\\text{(by
(\ref{pf.lem.sol.vander-hook.lap1.B.1}))}}}\right) \\
&  \ \ \ \ \ \ \ \ \ \ \left(
\begin{array}
[c]{c}%
\text{by Theorem \ref{thm.laplace.gen} \textbf{(a)} (applied to }B\text{,
}y_{i}^{n-j}\text{ and }k\\
\text{instead of }A\text{, }a_{i,j}\text{ and }p\text{)}%
\end{array}
\right) \\
&  =\sum_{q=1}^{n}\left(  -1\right)  ^{k+q}\left(  x_{k}+t\right)  ^{n-q}%
\det\left(  A_{\sim k,\sim q}\right)  .
\end{align*}
This proves Lemma \ref{lem.sol.vander-hook.lap1} \textbf{(b)}.

\textbf{(c)} For every $a\in\mathbb{K}$, $b\in\mathbb{K}$ and $m\in\mathbb{N}%
$, we have%
\begin{equation}
\left(  a+b\right)  ^{m}=\sum_{\ell=0}^{m}\dbinom{m}{\ell}a^{\ell}b^{m-\ell}.
\label{pf.lem.sol.vander-hook.lap1.c.binom}%
\end{equation}
(Indeed, this is precisely the equality (\ref{eq.rings.(a+b)**n}), with the
variables $n$ and $k$ renamed as $m$ and $\ell$.) Now, every $q\in\left\{
1,2,\ldots,n\right\}  $ satisfies%
\begin{equation}
\left(  x_{k}+t\right)  ^{n-q}-x_{k}^{n-q}=\sum_{\ell=1}^{n-q}\dbinom
{n-q}{\ell}t^{\ell}x_{k}^{n-q-\ell}
\label{pf.lem.sol.vander-hook.lap1.c.binom-used}%
\end{equation}
\footnote{\textit{Proof of (\ref{pf.lem.sol.vander-hook.lap1.c.binom-used}):}
Let $q\in\left\{  1,2,\ldots,n\right\}  $. Thus, $q\leq n$.
\par
Set $m=n-q$. Then, $m=n-q\in\mathbb{N}$ (since $q\leq n$). The equality
(\ref{pf.lem.sol.vander-hook.lap1.c.binom}) (applied to $a=t$ and $b=x_{k}$)
yields%
\begin{align*}
\left(  t+x_{k}\right)  ^{m}  &  =\sum_{\ell=0}^{m}\dbinom{m}{\ell}t^{\ell
}x_{k}^{m-\ell}=\underbrace{\dbinom{m}{0}}_{=1}\underbrace{t^{0}}%
_{=1}\underbrace{x_{k}^{m-0}}_{=x_{k}^{m}}+\sum_{\ell=1}^{m}\dbinom{m}{\ell
}t^{\ell}x_{k}^{m-\ell}\\
&  \ \ \ \ \ \ \ \ \ \ \left(  \text{here, we have split off the addend for
}\ell=0\text{ from the sum}\right) \\
&  =x_{k}^{m}+\sum_{\ell=1}^{m}\dbinom{m}{\ell}t^{\ell}x_{k}^{m-\ell}.
\end{align*}
Subtracting $x_{k}^{m}$ from both sides of this equality, we obtain%
\[
\left(  t+x_{k}\right)  ^{m}-x_{k}^{m}=\sum_{\ell=1}^{m}\dbinom{m}{\ell
}t^{\ell}x_{k}^{m-\ell}.
\]
Since $t+x_{k}=x_{k}+t$, this rewrites as
\[
\left(  x_{k}+t\right)  ^{m}-x_{k}^{m}=\sum_{\ell=1}^{m}\dbinom{m}{\ell
}t^{\ell}x_{k}^{m-\ell}.
\]
Since $m=n-q$, this rewrites as
\[
\left(  x_{k}+t\right)  ^{n-q}-x_{k}^{n-q}=\sum_{\ell=1}^{n-q}\dbinom
{n-q}{\ell}t^{\ell}x_{k}^{n-q-\ell}.
\]
This proves (\ref{pf.lem.sol.vander-hook.lap1.c.binom-used}).}.

But%
\begin{align*}
&  \underbrace{V\left(  x_{1},x_{2},\ldots,x_{k-1},x_{k}+t,x_{k+1}%
,x_{k+2},\ldots,x_{n}\right)  }_{\substack{=\sum_{q=1}^{n}\left(  -1\right)
^{k+q}\left(  x_{k}+t\right)  ^{n-q}\det\left(  A_{\sim k,\sim q}\right)
\\\text{(by Lemma \ref{lem.sol.vander-hook.lap1} \textbf{(b)})}}%
}-\underbrace{V\left(  x_{1},x_{2},\ldots,x_{n}\right)  }_{\substack{=\sum
_{q=1}^{n}\left(  -1\right)  ^{k+q}x_{k}^{n-q}\det\left(  A_{\sim k,\sim
q}\right)  \\\text{(by Lemma \ref{lem.sol.vander-hook.lap1} \textbf{(a)})}}}\\
&  =\sum_{q=1}^{n}\left(  -1\right)  ^{k+q}\left(  x_{k}+t\right)  ^{n-q}%
\det\left(  A_{\sim k,\sim q}\right)  -\sum_{q=1}^{n}\left(  -1\right)
^{k+q}x_{k}^{n-q}\det\left(  A_{\sim k,\sim q}\right) \\
&  =\sum_{q=1}^{n}\left(  -1\right)  ^{k+q}\underbrace{\left(  \left(
x_{k}+t\right)  ^{n-q}\det\left(  A_{\sim k,\sim q}\right)  -x_{k}^{n-q}%
\det\left(  A_{\sim k,\sim q}\right)  \right)  }_{=\left(  \left(
x_{k}+t\right)  ^{n-q}-x_{k}^{n-q}\right)  \det\left(  A_{\sim k,\sim
q}\right)  }\\
&  =\sum_{q=1}^{n}\left(  -1\right)  ^{k+q}\underbrace{\left(  \left(
x_{k}+t\right)  ^{n-q}-x_{k}^{n-q}\right)  }_{\substack{=\sum_{\ell=1}%
^{n-q}\dbinom{n-q}{\ell}t^{\ell}x_{k}^{n-q-\ell}\\\text{(by
(\ref{pf.lem.sol.vander-hook.lap1.c.binom-used}))}}}\det\left(  A_{\sim k,\sim
q}\right) \\
&  =\sum_{q=1}^{n}\left(  -1\right)  ^{k+q}\left(  \sum_{\ell=1}^{n-q}%
\dbinom{n-q}{\ell}t^{\ell}x_{k}^{n-q-\ell}\right)  \det\left(  A_{\sim k,\sim
q}\right) \\
&  =\sum_{q=1}^{n}\ \ \sum_{\ell=1}^{n-q}\underbrace{\left(  -1\right)
^{k+q}\dbinom{n-q}{\ell}t^{\ell}}_{=\dbinom{n-q}{\ell}t^{\ell}\left(
-1\right)  ^{k+q}}x_{k}^{n-q-\ell}\det\left(  A_{\sim k,\sim q}\right) \\
&  =\sum_{q=1}^{n}\ \ \sum_{\ell=1}^{n-q}\dbinom{n-q}{\ell}t^{\ell}\left(
-1\right)  ^{k+q}x_{k}^{n-q-\ell}\det\left(  A_{\sim k,\sim q}\right)  .
\end{align*}
This proves Lemma \ref{lem.sol.vander-hook.lap1} \textbf{(c)}.
\end{proof}

\begin{lemma}
\label{lem.sol.vander-hook.lap2}Let $n\in\mathbb{N}$. Let $x_{1},x_{2}%
,\ldots,x_{n}$ be $n$ elements of $\mathbb{K}$. Let $A$ be the $n\times
n$-matrix $\left(  x_{i}^{n-j}\right)  _{1\leq i\leq n,\ 1\leq j\leq n}$. Let
$q\in\left\{  1,2,\ldots,n\right\}  $. Let $\ell\in\left\{  1,2,\ldots
,n-q\right\}  $. Then,%
\[
\sum_{k=1}^{n}\left(  -1\right)  ^{k+q}x_{k}^{n-q-\ell+1}\det\left(  A_{\sim
k,\sim q}\right)  =\delta_{\ell,1}\det A.
\]

\end{lemma}

\begin{proof}
[Proof of Lemma \ref{lem.sol.vander-hook.lap2}.]From $\ell\in\left\{
1,2,\ldots,n-q\right\}  $, we obtain $\ell\geq1$ and $\ell\leq n-q$. Now,
$q+\underbrace{\ell}_{\leq n-q}-1\leq q+\left(  n-q\right)  -1=n-1\leq n$.
Also, $q\in\left\{  1,2,\ldots,n\right\}  $, so that $q\geq1$. Hence,
$\underbrace{q}_{\geq1}+\underbrace{\ell}_{\geq1}-1\geq1+1-1=1$. Combining
this with $q+\ell-1\leq n$, we obtain $q+\ell-1\in\left\{  1,2,\ldots
,n\right\}  $.

\begin{vershort}
We have $\delta_{q,q+\ell-1}=\delta_{\ell,1}$ (since $q=q+\ell-1$ holds if and
only if $\ell=1$).
\end{vershort}

\begin{verlong}
We have $\delta_{q,q+\ell-1}=\delta_{\ell,1}$\ \ \ \ \footnote{\textit{Proof.}
We are in one of the following two cases:
\par
\textit{Case 1:} We have $\ell=1$.
\par
\textit{Case 2:} We have $\ell\neq1$.
\par
Let us first consider Case 1. In this case, we have $\ell=1$. Thus,
$\delta_{\ell,1}=1$. But $q+\underbrace{\ell}_{=1}-1=q+1-1=q$, so that
$q=q+\ell-1$. Hence, $\delta_{q,q+\ell-1}=1$. Comparing this with
$\delta_{\ell,1}=1$, we obtain $\delta_{q,q+\ell-1}=\delta_{\ell,1}$. Hence,
$\delta_{q,q+\ell-1}=\delta_{\ell,1}$ is proven in Case 1.
\par
Let us first consider Case 2. In this case, we have $\ell\neq1$. Thus,
$\delta_{\ell,1}=0$. But $q+\ell-1\neq q$ (since $\left(  q+\ell-1\right)
-q=\ell-1\neq0$ (since $\ell\neq1$)). In other words, $q\neq q+\ell-1$. Hence,
$\delta_{q,q+\ell-1}=0$. Comparing this with $\delta_{\ell,1}=0$, we obtain
$\delta_{q,q+\ell-1}=\delta_{\ell,1}$. Hence, $\delta_{q,q+\ell-1}%
=\delta_{\ell,1}$ is proven in Case 2.
\par
We have now proven $\delta_{q,q+\ell-1}=\delta_{\ell,1}$ in each of the two
Cases 1 and 2. Since these two Cases cover all possibilities, this shows that
$\delta_{q,q+\ell-1}=\delta_{\ell,1}$ always holds. Qed.}.
\end{verlong}

We have $A=\left(  x_{i}^{n-j}\right)  _{1\leq i\leq n,\ 1\leq j\leq n}$ (by
the definition of $A$) and $q+\ell-1\in\left\{  1,2,\ldots,n\right\}  $.
Hence, Lemma \ref{lem.sol.vander-hook.lap0} (applied to $x_{i}^{n-j}$ and
$q+\ell-1$ instead of $a_{i,j}$ and $r$) yields%
\[
\sum_{p=1}^{n}\left(  -1\right)  ^{p+q}x_{p}^{n-\left(  q+\ell-1\right)  }%
\det\left(  A_{\sim p,\sim q}\right)  =\underbrace{\delta_{q,q+\ell-1}%
}_{=\delta_{\ell,1}}\det A=\delta_{\ell,1}\det A.
\]
Thus,%
\begin{align*}
\delta_{\ell,1}\det A  &  =\sum_{p=1}^{n}\left(  -1\right)  ^{p+q}%
\underbrace{x_{p}^{n-\left(  q+\ell-1\right)  }}_{\substack{=x_{p}%
^{n-q-\ell+1}\\\text{(since }n-\left(  q+\ell-1\right)  =n-q-\ell+1\text{)}%
}}\det\left(  A_{\sim p,\sim q}\right) \\
&  =\sum_{p=1}^{n}\left(  -1\right)  ^{p+q}x_{p}^{n-q-\ell+1}\det\left(
A_{\sim p,\sim q}\right)  =\sum_{k=1}^{n}\left(  -1\right)  ^{k+q}%
x_{k}^{n-q-\ell+1}\det\left(  A_{\sim k,\sim q}\right)
\end{align*}
(here, we have renamed the summation index $p$ as $k$). This proves Lemma
\ref{lem.sol.vander-hook.lap2}.
\end{proof}

\subsubsection{The solution}

\begin{vershort}
\begin{proof}
[Solution to Exercise \ref{exe.vander-hook}.]Let $A$ be the $n\times n$-matrix
$\left(  x_{i}^{n-j}\right)  _{1\leq i\leq n,\ 1\leq j\leq n}$.

We have%
\begin{align}
&  \sum_{k=1}^{n}x_{k}V\left(  x_{1},x_{2},\ldots,x_{k-1},x_{k}+t,x_{k+1}%
,x_{k+2},\ldots,x_{n}\right)  -\sum_{k=1}^{n}x_{k}V\left(  x_{1},x_{2}%
,\ldots,x_{n}\right) \nonumber\\
&  =\sum_{k=1}^{n}x_{k}\underbrace{\left(  V\left(  x_{1},x_{2},\ldots
,x_{k-1},x_{k}+t,x_{k+1},x_{k+2},\ldots,x_{n}\right)  -V\left(  x_{1}%
,x_{2},\ldots,x_{n}\right)  \right)  }_{\substack{=\sum_{q=1}^{n}%
\ \ \sum_{\ell=1}^{n-q}\dbinom{n-q}{\ell}t^{\ell}\left(  -1\right)
^{k+q}x_{k}^{n-q-\ell}\det\left(  A_{\sim k,\sim q}\right)  \\\text{(by Lemma
\ref{lem.sol.vander-hook.lap1} \textbf{(c)})}}}\nonumber\\
&  =\sum_{k=1}^{n}x_{k}\sum_{q=1}^{n}\ \ \sum_{\ell=1}^{n-q}\dbinom{n-q}{\ell
}t^{\ell}\left(  -1\right)  ^{k+q}x_{k}^{n-q-\ell}\det\left(  A_{\sim k,\sim
q}\right) \nonumber\\
&  =\sum_{q=1}^{n}\ \ \sum_{\ell=1}^{n-q}\dbinom{n-q}{\ell}t^{\ell}\sum
_{k=1}^{n}\underbrace{x_{k}\left(  -1\right)  ^{k+q}x_{k}^{n-q-\ell}%
}_{=\left(  -1\right)  ^{k+q}x_{k}^{n-q-\ell+1}}\det\left(  A_{\sim k,\sim
q}\right) \nonumber\\
&  =\sum_{q=1}^{n}\ \ \sum_{\ell=1}^{n-q}\dbinom{n-q}{\ell}t^{\ell
}\underbrace{\sum_{k=1}^{n}\left(  -1\right)  ^{k+q}x_{k}^{n-q-\ell+1}%
\det\left(  A_{\sim k,\sim q}\right)  }_{\substack{=\delta_{\ell,1}\det
A\\\text{(by Lemma \ref{lem.sol.vander-hook.lap2})}}}\nonumber\\
&  =\sum_{q=1}^{n}\ \ \sum_{\ell=1}^{n-q}\dbinom{n-q}{\ell}t^{\ell}%
\delta_{\ell,1}\det A=\sum_{q=1}^{n}\left(  \sum_{\ell=1}^{n-q}\dbinom
{n-q}{\ell}t^{\ell}\delta_{\ell,1}\right)  \det A\nonumber\\
&  =\sum_{r=0}^{n-1}\left(  \sum_{\ell=1}^{r}\dbinom{r}{\ell}t^{\ell}%
\delta_{\ell,1}\right)  \det A \label{sol.vander-hook.short.3}%
\end{align}
(here, we have substituted $r$ for $n-q$ in the outer sum).

But every $r\in\left\{  0,1,\ldots,n-1\right\}  $ satisfies%
\begin{equation}
\sum_{\ell=1}^{r}\dbinom{r}{\ell}t^{\ell}\delta_{\ell,1}=rt
\label{sol.vander-hook.short.suml}%
\end{equation}
\footnote{\textit{Proof of (\ref{sol.vander-hook.short.suml}):} Let
$r\in\left\{  0,1,\ldots,n-1\right\}  $. We must prove
(\ref{sol.vander-hook.short.suml}).
\par
We are in one of the following two cases:
\par
\textit{Case 1:} We have $r\neq0$.
\par
\textit{Case 2:} We have $r=0$.
\par
Let us consider Case 1 first. In this case, we have $r\neq0$. Hence, $r>0$.
Thus, $1\in\left\{  1,2,\ldots,r\right\}  $. Thus, we can split off the addend
for $\ell=1$ from the sum $\sum_{\ell=1}^{r}\dbinom{r}{\ell}t^{\ell}%
\delta_{\ell,1}$. As a result, we obtain
\[
\sum_{\ell=1}^{r}\dbinom{r}{\ell}t^{\ell}\delta_{\ell,1}=\underbrace{\dbinom
{r}{1}}_{=r}\underbrace{t^{1}}_{=t}\underbrace{\delta_{1,1}}_{=1}+\sum
_{\ell=2}^{r}\dbinom{r}{\ell}t^{\ell}\underbrace{\delta_{\ell,1}%
}_{\substack{=0\\\text{(since }\ell\neq1\text{)}}}=rt+\underbrace{\sum
_{\ell=2}^{r}\dbinom{r}{\ell}t^{\ell}0}_{=0}=rt.
\]
Thus, (\ref{sol.vander-hook.short.suml}) is proven in Case 1.
\par
Proving (\ref{sol.vander-hook.short.suml}) in Case 2 is straightforward and
left to the reader.
\par
We now have proven (\ref{sol.vander-hook.short.suml}) in each of the two Cases
1 and 2. Thus, (\ref{sol.vander-hook.short.suml}) always holds.}.

Now, (\ref{sol.vander-hook.short.3}) becomes%
\begin{align*}
&  \sum_{k=1}^{n}x_{k}V\left(  x_{1},x_{2},\ldots,x_{k-1},x_{k}+t,x_{k+1}%
,x_{k+2},\ldots,x_{n}\right)  -\sum_{k=1}^{n}x_{k}V\left(  x_{1},x_{2}%
,\ldots,x_{n}\right) \\
&  =\sum_{r=0}^{n-1}\underbrace{\left(  \sum_{\ell=1}^{r}\dbinom{r}{\ell
}t^{\ell}\delta_{\ell,1}\right)  }_{\substack{=rt\\\text{(by
(\ref{sol.vander-hook.short.suml}))}}}\det A=\sum_{r=0}^{n-1}rt\det A\\
&  =\underbrace{\left(  \sum_{r=0}^{n-1}r\right)  }_{\substack{=\dbinom{n}%
{2}\\\text{(by (\ref{sol.vander-hook.gauss}))}}}t\underbrace{\det
A}_{\substack{=V\left(  x_{1},x_{2},\ldots,x_{n}\right)  \\\text{(by
(\ref{eq.lem.sol.vander-hook.V=x.eq}))}}}=\dbinom{n}{2}tV\left(  x_{1}%
,x_{2},\ldots,x_{n}\right)  .
\end{align*}
If we add $\sum_{k=1}^{n}x_{k}V\left(  x_{1},x_{2},\ldots,x_{n}\right)  $ to
both sides of this equality, then we obtain%
\begin{align*}
&  \sum_{k=1}^{n}x_{k}V\left(  x_{1},x_{2},\ldots,x_{k-1},x_{k}+t,x_{k+1}%
,x_{k+2},\ldots,x_{n}\right) \\
&  =\dbinom{n}{2}tV\left(  x_{1},x_{2},\ldots,x_{n}\right)  +\sum_{k=1}%
^{n}x_{k}V\left(  x_{1},x_{2},\ldots,x_{n}\right) \\
&  =\left(  \dbinom{n}{2}t+\sum_{k=1}^{n}x_{k}\right)  V\left(  x_{1}%
,x_{2},\ldots,x_{n}\right)  .
\end{align*}
This solves Exercise \ref{exe.vander-hook}.
\end{proof}
\end{vershort}

\begin{verlong}
\begin{proof}
[Solution to Exercise \ref{exe.vander-hook}.]Let $A$ be the $n\times n$-matrix
$\left(  x_{i}^{n-j}\right)  _{1\leq i\leq n,\ 1\leq j\leq n}$.

We have%
\begin{align}
&  \sum_{k=1}^{n}x_{k}V\left(  x_{1},x_{2},\ldots,x_{k-1},x_{k}+t,x_{k+1}%
,x_{k+2},\ldots,x_{n}\right) \nonumber\\
&  \ \ \ \ \ \ \ \ \ \ -\sum_{k=1}^{n}x_{k}V\left(  x_{1},x_{2},\ldots
,x_{n}\right) \nonumber\\
&  =\sum_{k=1}^{n}x_{k}\underbrace{\left(  V\left(  x_{1},x_{2},\ldots
,x_{k-1},x_{k}+t,x_{k+1},x_{k+2},\ldots,x_{n}\right)  -V\left(  x_{1}%
,x_{2},\ldots,x_{n}\right)  \right)  }_{\substack{=\sum_{q=1}^{n}%
\ \ \sum_{\ell=1}^{n-q}\dbinom{n-q}{\ell}t^{\ell}\left(  -1\right)
^{k+q}x_{k}^{n-q-\ell}\det\left(  A_{\sim k,\sim q}\right)  \\\text{(by Lemma
\ref{lem.sol.vander-hook.lap1} \textbf{(c)})}}}\nonumber\\
&  =\sum_{k=1}^{n}\underbrace{x_{k}\sum_{q=1}^{n}\ \ \sum_{\ell=1}%
^{n-q}\dbinom{n-q}{\ell}t^{\ell}\left(  -1\right)  ^{k+q}x_{k}^{n-q-\ell}%
\det\left(  A_{\sim k,\sim q}\right)  }_{=\sum_{q=1}^{n}\sum_{\ell=1}%
^{n-q}x_{k}\dbinom{n-q}{\ell}t^{\ell}\left(  -1\right)  ^{k+q}x_{k}^{n-q-\ell
}\det\left(  A_{\sim k,\sim q}\right)  }\nonumber\\
&  =\underbrace{\sum_{k=1}^{n}\ \ \sum_{q=1}^{n}\ \ \sum_{\ell=1}^{n-q}%
}_{=\sum_{q=1}^{n}\ \ \sum_{\ell=1}^{n-q}\ \ \sum_{k=1}^{n}}\underbrace{x_{k}%
\dbinom{n-q}{\ell}t^{\ell}\left(  -1\right)  ^{k+q}}_{=\dbinom{n-q}{\ell
}t^{\ell}\left(  -1\right)  ^{k+q}x_{k}}x_{k}^{n-q-\ell}\det\left(  A_{\sim
k,\sim q}\right) \nonumber\\
&  =\sum_{q=1}^{n}\ \ \sum_{\ell=1}^{n-q}\ \ \sum_{k=1}^{n}\dbinom{n-q}{\ell
}t^{\ell}\left(  -1\right)  ^{k+q}\underbrace{x_{k}x_{k}^{n-q-\ell}}%
_{=x_{k}^{n-q-\ell+1}}\det\left(  A_{\sim k,\sim q}\right) \nonumber\\
&  =\sum_{q=1}^{n}\ \ \sum_{\ell=1}^{n-q}\ \ \underbrace{\sum_{k=1}^{n}%
\dbinom{n-q}{\ell}t^{\ell}\left(  -1\right)  ^{k+q}x_{k}^{n-q-\ell+1}%
\det\left(  A_{\sim k,\sim q}\right)  }_{=\dbinom{n-q}{\ell}t^{\ell}\sum
_{k=1}^{n}\left(  -1\right)  ^{k+q}x_{k}^{n-q-\ell+1}\det\left(  A_{\sim
k,\sim q}\right)  }\nonumber\\
&  =\sum_{q=1}^{n}\ \ \sum_{\ell=1}^{n-q}\dbinom{n-q}{\ell}t^{\ell
}\underbrace{\sum_{k=1}^{n}\left(  -1\right)  ^{k+q}x_{k}^{n-q-\ell+1}%
\det\left(  A_{\sim k,\sim q}\right)  }_{\substack{=\delta_{\ell,1}\det
A\\\text{(by Lemma \ref{lem.sol.vander-hook.lap2})}}}\nonumber\\
&  =\sum_{q=1}^{n}\ \ \sum_{\ell=1}^{n-q}\dbinom{n-q}{\ell}t^{\ell}%
\delta_{\ell,1}\det A=\sum_{q=1}^{n}\left(  \sum_{\ell=1}^{n-q}\dbinom
{n-q}{\ell}t^{\ell}\delta_{\ell,1}\right)  \det A\nonumber\\
&  =\sum_{r=0}^{n-1}\left(  \sum_{\ell=1}^{r}\dbinom{r}{\ell}t^{\ell}%
\delta_{\ell,1}\right)  \det A \label{sol.vander-hook.3}%
\end{align}
(here, we have substituted $r$ for $n-q$ in the outer sum).

But every $r\in\left\{  0,1,\ldots,n-1\right\}  $ satisfies%
\begin{equation}
\sum_{\ell=1}^{r}\dbinom{r}{\ell}t^{\ell}\delta_{\ell,1}=rt
\label{sol.vander-hook.suml}%
\end{equation}
\footnote{\textit{Proof of (\ref{sol.vander-hook.suml}):} Let $r\in\left\{
0,1,\ldots,n-1\right\}  $. We must prove (\ref{sol.vander-hook.suml}).
\par
We are in one of the following two cases:
\par
\textit{Case 1:} We have $r\neq0$.
\par
\textit{Case 2:} We have $r=0$.
\par
Let us consider Case 1 first. In this case, we have $r\neq0$. But $r\geq0$
(since $r\in\left\{  0,1,\ldots,n-1\right\}  $). Hence, $r>0$ (since $r\neq
0$). Thus, $1\in\left\{  1,2,\ldots,r\right\}  $. Thus, we can split off the
addend for $\ell=1$ from the sum $\sum_{\ell=1}^{r}\dbinom{r}{\ell}t^{\ell
}\delta_{\ell,1}$. As a result, we obtain
\[
\sum_{\ell=1}^{r}\dbinom{r}{\ell}t^{\ell}\delta_{\ell,1}=\underbrace{\dbinom
{r}{1}}_{=r}\underbrace{t^{1}}_{=t}\underbrace{\delta_{1,1}}%
_{\substack{=1\\\text{(since }1=1\text{)}}}+\sum_{\ell=2}^{r}\dbinom{r}{\ell
}t^{\ell}\underbrace{\delta_{\ell,1}}_{\substack{=0\\\text{(since }\ell
\neq1\\\text{(since }\ell\geq2>1\text{))}}}=rt+\underbrace{\sum_{\ell=2}%
^{r}\dbinom{r}{\ell}t^{\ell}0}_{=0}=rt.
\]
Thus, (\ref{sol.vander-hook.suml}) is proven in Case 1.
\par
Let us now consider Case 2. In this case, we have $r=0$. Hence,%
\[
\sum_{\ell=1}^{r}\dbinom{r}{\ell}t^{\ell}\delta_{\ell,1}=\sum_{\ell=1}%
^{0}\dbinom{0}{\ell}t^{\ell}\delta_{\ell,1}=\left(  \text{empty sum}\right)
=0=rt
\]
(since $\underbrace{r}_{=0}t=0t=0$). Thus, (\ref{sol.vander-hook.suml}) is
proven in Case 2.
\par
We now have proven (\ref{sol.vander-hook.suml}) in each of the two Cases 1 and
2. Since these two Cases cover all possibilities, this shows that
(\ref{sol.vander-hook.suml}) always holds. Qed.}.

Now, (\ref{sol.vander-hook.3}) becomes%
\begin{align*}
&  \sum_{k=1}^{n}x_{k}V\left(  x_{1},x_{2},\ldots,x_{k-1},x_{k}+t,x_{k+1}%
,x_{k+2},\ldots,x_{n}\right) \\
&  \ \ \ \ \ \ \ \ \ \ -\sum_{k=1}^{n}x_{k}V\left(  x_{1},x_{2},\ldots
,x_{n}\right) \\
&  =\sum_{r=0}^{n-1}\underbrace{\left(  \sum_{\ell=1}^{r}\dbinom{r}{\ell
}t^{\ell}\delta_{\ell,1}\right)  }_{\substack{=rt\\\text{(by
(\ref{sol.vander-hook.suml}))}}}\det A=\sum_{r=0}^{n-1}rt\det A\\
&  =\underbrace{\left(  \sum_{r=0}^{n-1}r\right)  }_{\substack{=\dbinom{n}%
{2}\\\text{(by (\ref{sol.vander-hook.gauss}))}}}t\underbrace{\det
A}_{\substack{=V\left(  x_{1},x_{2},\ldots,x_{n}\right)  \\\text{(by
(\ref{eq.lem.sol.vander-hook.V=x.eq}))}}}\\
&  =\dbinom{n}{2}tV\left(  x_{1},x_{2},\ldots,x_{n}\right)  .
\end{align*}
If we add $\sum_{k=1}^{n}x_{k}V\left(  x_{1},x_{2},\ldots,x_{n}\right)  $ to
both sides of this equality, then we obtain%
\begin{align*}
&  \sum_{k=1}^{n}x_{k}V\left(  x_{1},x_{2},\ldots,x_{k-1},x_{k}+t,x_{k+1}%
,x_{k+2},\ldots,x_{n}\right) \\
&  =\dbinom{n}{2}tV\left(  x_{1},x_{2},\ldots,x_{n}\right)  +\sum_{k=1}%
^{n}x_{k}V\left(  x_{1},x_{2},\ldots,x_{n}\right) \\
&  =\left(  \dbinom{n}{2}t+\sum_{k=1}^{n}x_{k}\right)  V\left(  x_{1}%
,x_{2},\ldots,x_{n}\right)  .
\end{align*}
This solves Exercise \ref{exe.vander-hook}.
\end{proof}
\end{verlong}

\subsubsection{Addendum: a simpler variant}

Exercise \ref{exe.vander-hook} is now solved, but let me discuss one further
fact, which is a variation on it. Namely, the following holds:

\begin{proposition}
\label{prop.sol.vander-hook.variant}Let $n\in\mathbb{N}$. Let $x_{1}%
,x_{2},\ldots,x_{n}$ be $n$ elements of $\mathbb{K}$. Let $t\in\mathbb{K}$.
Then,%
\[
\sum_{k=1}^{n}V\left(  x_{1},x_{2},\ldots,x_{k-1},x_{k}+t,x_{k+1}%
,x_{k+2},\ldots,x_{n}\right)  =nV\left(  x_{1},x_{2},\ldots,x_{n}\right)  .
\]

\end{proposition}

We can prove this similarly to how we solved Exercise \ref{exe.vander-hook}.
Instead of Lemma \ref{lem.sol.vander-hook.lap2}, we use the following variant
of this lemma:

\begin{lemma}
\label{lem.sol.vander-hook.lap2.variant}Let $n\in\mathbb{N}$. Let $x_{1}%
,x_{2},\ldots,x_{n}$ be $n$ elements of $\mathbb{K}$. Let $A$ be the $n\times
n$-matrix $\left(  x_{i}^{n-j}\right)  _{1\leq i\leq n,\ 1\leq j\leq n}$. Let
$q\in\left\{  1,2,\ldots,n\right\}  $. Let $\ell\in\left\{  1,2,\ldots
,n-q\right\}  $. Then,%
\[
\sum_{k=1}^{n}\left(  -1\right)  ^{k+q}x_{k}^{n-q-\ell}\det\left(  A_{\sim
k,\sim q}\right)  =0.
\]

\end{lemma}

\begin{vershort}
\begin{proof}
[Proof of Lemma \ref{lem.sol.vander-hook.lap2.variant}.]Left to the reader.
(Very similar to the above proof of Lemma \ref{lem.sol.vander-hook.lap2}.)
\end{proof}
\end{vershort}

\begin{verlong}
\begin{proof}
[Proof of Lemma \ref{lem.sol.vander-hook.lap2.variant}.]From $\ell\in\left\{
1,2,\ldots,n-q\right\}  $, we obtain $\ell\geq1$ and $\ell\leq n-q$. Now,
$q+\underbrace{\ell}_{\leq n-q}\leq q+\left(  n-q\right)  =n$. Also,
$q\in\left\{  1,2,\ldots,n\right\}  $, so that $q\geq1$. Hence,
$\underbrace{q}_{\geq1}+\underbrace{\ell}_{\geq1}\geq1+1\geq1$. Combining this
with $q+\ell\leq n$, we obtain $q+\ell\in\left\{  1,2,\ldots,n\right\}  $.

From $q+\underbrace{\ell}_{\geq1>0}>q$, we obtain $q+\ell\neq q$. Hence,
$q\neq q+\ell$ and thus $\delta_{q,q+\ell}=0$.

We have $A=\left(  x_{i}^{n-j}\right)  _{1\leq i\leq n,\ 1\leq j\leq n}$ (by
the definition of $A$) and $q+\ell\in\left\{  1,2,\ldots,n\right\}  $. Hence,
Lemma \ref{lem.sol.vander-hook.lap0} (applied to $x_{i}^{n-j}$ and $q+\ell$
instead of $a_{i,j}$ and $r$) yields%
\[
\sum_{p=1}^{n}\left(  -1\right)  ^{p+q}x_{p}^{n-\left(  q+\ell\right)  }%
\det\left(  A_{\sim p,\sim q}\right)  =\underbrace{\delta_{q,q+\ell}}_{=0}\det
A=0\det A=0.
\]
Thus,%
\begin{align*}
0  &  =\sum_{p=1}^{n}\left(  -1\right)  ^{p+q}\underbrace{x_{p}^{n-\left(
q+\ell\right)  }}_{\substack{=x_{p}^{n-q-\ell}\\\text{(since }n-\left(
q+\ell\right)  =n-q-\ell\text{)}}}\det\left(  A_{\sim p,\sim q}\right) \\
&  =\sum_{p=1}^{n}\left(  -1\right)  ^{p+q}x_{p}^{n-q-\ell}\det\left(  A_{\sim
p,\sim q}\right)  =\sum_{k=1}^{n}\left(  -1\right)  ^{k+q}x_{k}^{n-q-\ell}%
\det\left(  A_{\sim k,\sim q}\right)
\end{align*}
(here, we have renamed the summation index $p$ as $k$). Lemma
\ref{lem.sol.vander-hook.lap2.variant} is therefore proven.
\end{proof}
\end{verlong}

\begin{vershort}
\begin{proof}
[Proof of Proposition \ref{prop.sol.vander-hook.variant}.]This proof is
similar to our solution of Exercise \ref{exe.vander-hook}, but a lot simpler.
Again, the reader can fill in the details.
\end{proof}
\end{vershort}

\begin{verlong}
\begin{proof}
[Proof of Proposition \ref{prop.sol.vander-hook.variant}.]Let $A$ be the
$n\times n$-matrix $\left(  x_{i}^{n-j}\right)  _{1\leq i\leq n,\ 1\leq j\leq
n}$. We have%
\begin{align}
&  \sum_{k=1}^{n}V\left(  x_{1},x_{2},\ldots,x_{k-1},x_{k}+t,x_{k+1}%
,x_{k+2},\ldots,x_{n}\right)  -\sum_{k=1}^{n}V\left(  x_{1},x_{2},\ldots
,x_{n}\right) \nonumber\\
&  =\sum_{k=1}^{n}\underbrace{\left(  V\left(  x_{1},x_{2},\ldots
,x_{k-1},x_{k}+t,x_{k+1},x_{k+2},\ldots,x_{n}\right)  -V\left(  x_{1}%
,x_{2},\ldots,x_{n}\right)  \right)  }_{\substack{=\sum_{q=1}^{n}%
\ \ \sum_{\ell=1}^{n-q}\dbinom{n-q}{\ell}t^{\ell}\left(  -1\right)
^{k+q}x_{k}^{n-q-\ell}\det\left(  A_{\sim k,\sim q}\right)  \\\text{(by Lemma
\ref{lem.sol.vander-hook.lap1} \textbf{(c)})}}}\nonumber\\
&  =\underbrace{\sum_{k=1}^{n}\ \ \sum_{q=1}^{n}\ \ \sum_{\ell=1}^{n-q}%
}_{=\sum_{q=1}^{n}\ \ \sum_{\ell=1}^{n-q}\ \ \sum_{k=1}^{n}}\dbinom{n-q}{\ell
}t^{\ell}\left(  -1\right)  ^{k+q}x_{k}^{n-q-\ell}\det\left(  A_{\sim k,\sim
q}\right) \nonumber\\
&  =\sum_{q=1}^{n}\ \ \sum_{\ell=1}^{n-q}\ \ \underbrace{\sum_{k=1}^{n}%
\dbinom{n-q}{\ell}t^{\ell}\left(  -1\right)  ^{k+q}x_{k}^{n-q-\ell}\det\left(
A_{\sim k,\sim q}\right)  }_{=\dbinom{n-q}{\ell}t^{\ell}\sum_{k=1}^{n}\left(
-1\right)  ^{k+q}x_{k}^{n-q-\ell}\det\left(  A_{\sim k,\sim q}\right)
}\nonumber\\
&  =\sum_{q=1}^{n}\ \ \sum_{\ell=1}^{n-q}\dbinom{n-q}{\ell}t^{\ell
}\underbrace{\sum_{k=1}^{n}\left(  -1\right)  ^{k+q}x_{k}^{n-q-\ell}%
\det\left(  A_{\sim k,\sim q}\right)  }_{\substack{=0\\\text{(by Lemma
\ref{lem.sol.vander-hook.lap2.variant})}}}\nonumber\\
&  =\sum_{q=1}^{n}\ \ \sum_{\ell=1}^{n-q}\dbinom{n-q}{\ell}t^{\ell
}0=0.\nonumber
\end{align}
In other words,%
\begin{align*}
&  \sum_{k=1}^{n}V\left(  x_{1},x_{2},\ldots,x_{k-1},x_{k}+t,x_{k+1}%
,x_{k+2},\ldots,x_{n}\right) \\
&  =\sum_{k=1}^{n}V\left(  x_{1},x_{2},\ldots,x_{n}\right)  =nV\left(
x_{1},x_{2},\ldots,x_{n}\right)  .
\end{align*}
This proves Proposition \ref{prop.sol.vander-hook.variant}.
\end{proof}
\end{verlong}

\subsubsection{Addendum: another sum of Vandermonde determinants}

Here is one more result similar to Exercise \ref{exe.vander-hook} and
Proposition \ref{prop.sol.vander-hook.variant}:

\begin{proposition}
\label{prop.sol.vander-hook.N}Let $n\in\mathbb{N}$. Let $x_{1},x_{2}%
,\ldots,x_{n}$ be $n$ elements of $\mathbb{K}$. Let $m\in\left\{
0,1,\ldots,n-1\right\}  $. Let $t\in\mathbb{K}$. Then,%
\[
\sum_{k=1}^{n}x_{k}^{m}V\left(  x_{1},x_{2},\ldots,x_{k-1},t,x_{k+1}%
,x_{k+2},\ldots,x_{n}\right)  =t^{m}V\left(  x_{1},x_{2},\ldots,x_{n}\right)
.
\]

\end{proposition}

We have already built all the tools necessary for the proof of this
proposition. We just need to repurpose a few of them:

\begin{lemma}
\label{lem.sol.vander-hook.N.lap1b}Let $n\in\mathbb{N}$. Let $x_{1}%
,x_{2},\ldots,x_{n}$ be $n$ elements of $\mathbb{K}$. Let $t\in\mathbb{K}$.
Let $A$ be the $n\times n$-matrix $\left(  x_{i}^{n-j}\right)  _{1\leq i\leq
n,\ 1\leq j\leq n}$. Let $k\in\left\{  1,2,\ldots,n\right\}  $. Then,%
\[
V\left(  x_{1},x_{2},\ldots,x_{k-1},t,x_{k+1},x_{k+2},\ldots,x_{n}\right)
=\sum_{q=1}^{n}\left(  -1\right)  ^{k+q}t^{n-q}\det\left(  A_{\sim k,\sim
q}\right)  .
\]

\end{lemma}

\begin{proof}
[Proof of Lemma \ref{lem.sol.vander-hook.N.lap1b}.]Lemma
\ref{lem.sol.vander-hook.lap1} \textbf{(b)} (applied to $t-x_{k}$ instead of
$t$) yields%
\begin{align*}
&  V\left(  x_{1},x_{2},\ldots,x_{k-1},x_{k}+\left(  t-x_{k}\right)
,x_{k+1},x_{k+2},\ldots,x_{n}\right) \\
&  =\sum_{q=1}^{n}\left(  -1\right)  ^{k+q}\left(  x_{k}+\left(
t-x_{k}\right)  \right)  ^{n-q}\det\left(  A_{\sim k,\sim q}\right)  .
\end{align*}
Since $x_{k}+\left(  t-x_{k}\right)  =t$, this rewrites as
\[
V\left(  x_{1},x_{2},\ldots,x_{k-1},t,x_{k+1},x_{k+2},\ldots,x_{n}\right)
=\sum_{q=1}^{n}\left(  -1\right)  ^{k+q}t^{n-q}\det\left(  A_{\sim k,\sim
q}\right)  .
\]
This proves Lemma \ref{lem.sol.vander-hook.N.lap1b}.
\end{proof}

\begin{lemma}
\label{lem.sol.vander-hook.N.lap2}Let $n\in\mathbb{N}$. Let $x_{1}%
,x_{2},\ldots,x_{n}$ be $n$ elements of $\mathbb{K}$. Let $A$ be the $n\times
n$-matrix $\left(  x_{i}^{n-j}\right)  _{1\leq i\leq n,\ 1\leq j\leq n}$. Let
$q\in\left\{  1,2,\ldots,n\right\}  $ and $m\in\left\{  0,1,\ldots
,n-1\right\}  $. Then,%
\[
\sum_{k=1}^{n}x_{k}^{m}\left(  -1\right)  ^{k+q}\det\left(  A_{\sim k,\sim
q}\right)  =\delta_{q,n-m}\det A.
\]

\end{lemma}

\begin{proof}
[Proof of Lemma \ref{lem.sol.vander-hook.N.lap2}.]From $q\in\left\{
1,2,\ldots,n\right\}  $, we obtain $1\leq q\leq n$. Hence, $1\leq n$, so that
$n\geq1$ and therefore $n\in\left\{  1,2,\ldots,n\right\}  $.

We have $A=\left(  x_{i}^{n-j}\right)  _{1\leq i\leq n,\ 1\leq j\leq n}$ (by
the definition of $A$) and $n-m\in\left\{  1,2,\ldots,n\right\}  $ (since
$m\in\left\{  0,1,\ldots,n-1\right\}  $). Hence, Lemma
\ref{lem.sol.vander-hook.lap0} (applied to $x_{i}^{n-j}$ and $n-m$ instead of
$a_{i,j}$ and $r$) yields%
\[
\sum_{p=1}^{n}\left(  -1\right)  ^{p+q}x_{p}^{n-\left(  n-m\right)  }%
\det\left(  A_{\sim p,\sim q}\right)  =\delta_{q,n-m}\det A.
\]
Thus,%
\begin{align*}
\delta_{q,n-m}\det A  &  =\sum_{p=1}^{n}\left(  -1\right)  ^{p+q}%
\underbrace{x_{p}^{n-\left(  n-m\right)  }}_{\substack{=x_{p}^{m}%
\\\text{(since }n-\left(  n-m\right)  =m\text{)}}}\det\left(  A_{\sim p,\sim
q}\right) \\
&  =\sum_{p=1}^{n}\underbrace{\left(  -1\right)  ^{p+q}x_{p}^{m}}_{=x_{p}%
^{m}\left(  -1\right)  ^{p+q}}\det\left(  A_{\sim p,\sim q}\right) \\
&  =\sum_{p=1}^{n}x_{p}^{m}\left(  -1\right)  ^{p+q}\det\left(  A_{\sim p,\sim
q}\right)  =\sum_{k=1}^{n}x_{k}^{m}\left(  -1\right)  ^{k+q}\det\left(
A_{\sim k,\sim q}\right)
\end{align*}
(here, we have renamed the summation index $p$ as $k$). This proves Lemma
\ref{lem.sol.vander-hook.N.lap2}.
\end{proof}

\begin{proof}
[Proof of Proposition \ref{prop.sol.vander-hook.N}.]Let $A$ be the $n\times
n$-matrix $\left(  x_{i}^{n-j}\right)  _{1\leq i\leq n,\ 1\leq j\leq n}$.

From $m\in\left\{  0,1,\ldots,n-1\right\}  $, we obtain $n-m\in\left\{
1,2,\ldots,n\right\}  $. But%
\begin{align*}
&  \sum_{k=1}^{n}x_{k}^{m}\underbrace{V\left(  x_{1},x_{2},\ldots
,x_{k-1},t,x_{k+1},x_{k+2},\ldots,x_{n}\right)  }_{\substack{=\sum_{q=1}%
^{n}\left(  -1\right)  ^{k+q}t^{n-q}\det\left(  A_{\sim k,\sim q}\right)
\\\text{(by Lemma \ref{lem.sol.vander-hook.N.lap1b})}}}\\
&  =\sum_{k=1}^{n}\underbrace{x_{k}^{m}\sum_{q=1}^{n}\left(  -1\right)
^{k+q}t^{n-q}\det\left(  A_{\sim k,\sim q}\right)  }_{=\sum_{q=1}^{n}x_{k}%
^{m}\left(  -1\right)  ^{k+q}t^{n-q}\det\left(  A_{\sim k,\sim q}\right)  }\\
&  =\underbrace{\sum_{k=1}^{n}\ \ \sum_{q=1}^{n}}_{=\sum_{q=1}^{n}%
\ \ \sum_{k=1}^{n}}x_{k}^{m}\left(  -1\right)  ^{k+q}t^{n-q}\det\left(
A_{\sim k,\sim q}\right) \\
&  =\sum_{q=1}^{n}\ \ \underbrace{\sum_{k=1}^{n}x_{k}^{m}\left(  -1\right)
^{k+q}t^{n-q}\det\left(  A_{\sim k,\sim q}\right)  }_{=t^{n-q}\sum_{k=1}%
^{n}x_{k}^{m}\left(  -1\right)  ^{k+q}\det\left(  A_{\sim k,\sim q}\right)
}\\
&  =\underbrace{\sum_{q=1}^{n}}_{=\sum_{q\in\left\{  1,2,\ldots,n\right\}  }%
}t^{n-q}\underbrace{\sum_{k=1}^{n}x_{k}^{m}\left(  -1\right)  ^{k+q}%
\det\left(  A_{\sim k,\sim q}\right)  }_{\substack{=\delta_{q,n-m}\det
A\\\text{(by Lemma \ref{lem.sol.vander-hook.N.lap2})}}}=\sum_{q\in\left\{
1,2,\ldots,n\right\}  }t^{n-q}\delta_{q,n-m}\det A\\
&  =\sum_{\substack{q\in\left\{  1,2,\ldots,n\right\}  ;\\q\neq n-m}%
}t^{n-q}\underbrace{\delta_{q,n-m}}_{\substack{=0\\\text{(since }q\neq
n-m\text{)}}}\det A+\underbrace{t^{n-\left(  n-m\right)  }}_{\substack{=t^{m}%
\\\text{(since }n-\left(  n-m\right)  =m\text{)}}}\underbrace{\delta
_{n-m,n-m}}_{\substack{=1\\\text{(since }n-m=n-m\text{)}}}\det A\\
&  \ \ \ \ \ \ \ \ \ \ \left(
\begin{array}
[c]{c}%
\text{here, we have split off the addend for }q=n-m\text{ from the sum}\\
\text{(since }n-m\in\left\{  1,2,\ldots,n\right\}  \text{)}%
\end{array}
\right) \\
&  =\underbrace{\sum_{\substack{q\in\left\{  1,2,\ldots,n\right\}  ;\\q\neq
n-m}}t^{n-q}0\det A}_{=0}+t^{m}\det A=t^{m}\underbrace{\det A}%
_{\substack{=V\left(  x_{1},x_{2},\ldots,x_{n}\right)  \\\text{(by
(\ref{eq.lem.sol.vander-hook.V=x.eq}))}}}\\
&  =t^{m}V\left(  x_{1},x_{2},\ldots,x_{n}\right)  .
\end{align*}
This proves Proposition \ref{prop.sol.vander-hook.N}.
\end{proof}

\subsubsection{Addendum: analogues involving products of all but one $x_{j}$}

Let us finally prove a much more complicated analogue of Exercise
\ref{exe.vander-hook} and Proposition \ref{prop.sol.vander-hook.variant}. We
shall use the following notations:

\begin{definition}
\label{def.sol.vander-hook.elsyms}Let $n\in\mathbb{N}$. Let $\left[  n\right]
$ denote the set $\left\{  1,2,\ldots,n\right\}  $. As usual, let
$\mathcal{P}\left(  \left[  n\right]  \right)  $ denote the powerset of
$\left[  n\right]  $. Let $x_{1},x_{2},\ldots,x_{n}$ be $n$ elements of
$\mathbb{K}$.

\textbf{(a)} For every $j\in\mathbb{N}$, define an element $e_{j}\left(
x_{1},x_{2},\ldots,x_{n}\right)  \in\mathbb{K}$ by
\[
e_{j}\left(  x_{1},x_{2},\ldots,x_{n}\right)  =\sum_{\substack{I\subseteq
\left[  n\right]  ;\\\left\vert I\right\vert =j}}\ \ \prod_{i\in I}x_{i}.
\]

(Here, as usual, the summation sign $\sum_{\substack{I\subseteq\left[
n\right]  ;\\\left\vert I\right\vert =j}}$ means $\sum_{\substack{I\in
\mathcal{P}\left(  \left[  n\right]  \right)  ;\\\left\vert I\right\vert =j}}$.)

\textbf{(b)} For every $t\in\mathbb{K}$, define an element $z_{t}\left(
x_{1},x_{2},\ldots,x_{n}\right)  \in\mathbb{K}$ by%
\[
z_{t}\left(  x_{1},x_{2},\ldots,x_{n}\right)  =\sum_{j=0}^{n-1}e_{n-1-j}%
\left(  x_{1},x_{2},\ldots,x_{n}\right)  t^{j}.
\]

\end{definition}

\begin{proposition}
\label{prop.sol.vander-hook.variant-1}Let $n\in\mathbb{N}$. Let $x_{1}%
,x_{2},\ldots,x_{n}$ be $n$ elements of $\mathbb{K}$. Let $t\in\mathbb{K}$.
For each $i\in\left\{  1,2,\ldots,n\right\}  $, set $y_{i}=\prod
_{\substack{j\in\left\{  1,2,\ldots,n\right\}  ;\\j\neq i}}x_{j}$. Then,%
\begin{align*}
&  \sum_{k=1}^{n}y_{k}V\left(  x_{1},x_{2},\ldots,x_{k-1},x_{k}+t,x_{k+1}%
,x_{k+2},\ldots,x_{n}\right) \\
&  =z_{-t}\left(  x_{1},x_{2},\ldots,x_{n}\right)  \cdot V\left(  x_{1}%
,x_{2},\ldots,x_{n}\right)  .
\end{align*}

\end{proposition}

Before we start proving Proposition \ref{prop.sol.vander-hook.variant-1}, let
us explore a few properties of the elements defined in Definition
\ref{def.sol.vander-hook.elsyms}:

\begin{proposition}
\label{prop.sol.vander-hook.viete}Let $n\in\mathbb{N}$. Let $x_{1}%
,x_{2},\ldots,x_{n}$ be $n$ elements of $\mathbb{K}$. Let $t\in\mathbb{K}$.

\textbf{(a)} We have%
\[
\prod_{i=1}^{n}\left(  x_{i}+t\right)  =\sum_{j=0}^{n}e_{n-j}\left(
x_{1},x_{2},\ldots,x_{n}\right)  t^{j}.
\]

\textbf{(b)} We have%
\[
e_{n}\left(  x_{1},x_{2},\ldots,x_{n}\right)  =\prod_{i=1}^{n}x_{i}.
\]

\textbf{(c)} We have%
\[
t\cdot z_{t}\left(  x_{1},x_{2},\ldots,x_{n}\right)  =\prod_{i=1}^{n}\left(
x_{i}+t\right)  -\prod_{i=1}^{n}x_{i}.
\]

\textbf{(d)} Assume that the element $t$ of $\mathbb{K}$ is invertible. Then,%
\[
z_{t}\left(  x_{1},x_{2},\ldots,x_{n}\right)  =\dfrac{\prod_{i=1}^{n}\left(
x_{i}+t\right)  -\prod_{i=1}^{n}x_{i}}{t}.
\]

\end{proposition}

Proposition \ref{prop.sol.vander-hook.viete} \textbf{(a)} is one of several
interconnected results known as \textit{Vieta's formulas}. Proposition
\ref{prop.sol.vander-hook.viete} \textbf{(d)} gives an alternative description
of $z_{t}\left(  x_{1},x_{2},\ldots,x_{n}\right)  $ in the case when $t$ is invertible.

\begin{proof}
[Proof of Proposition \ref{prop.sol.vander-hook.viete}.]We shall use the
notations $\left[  n\right]  $, $\mathcal{P}\left(  \left[  n\right]  \right)
$ and $\sum_{\substack{I\subseteq\left[  n\right]  ;\\\left\vert I\right\vert
=j}}$ as in Definition \ref{def.sol.vander-hook.elsyms}.

\begin{vershort}
We have $\left\vert \left[  n\right]  \right\vert =n$. Hence, every subset $I$
of $\left[  n\right]  $ satisfies $\left\vert I\right\vert \in\left\{
0,1,\ldots,n\right\}  $.
\end{vershort}

\begin{verlong}
Every subset $I$ of $\left[  n\right]  $ satisfies $\left\vert I\right\vert
\in\left\{  0,1,\ldots,n\right\}  $\ \ \ \ \footnote{\textit{Proof.} Let $I$
be a subset of $\left[  n\right]  $. Then, $I\subseteq\left[  n\right]  $, so
that $\left\vert I\right\vert \leq\left\vert \left[  n\right]  \right\vert
=n$. In particular, $\left\vert I\right\vert \in\mathbb{N}$. Thus, $\left\vert
I\right\vert \in\left\{  0,1,\ldots,n\right\}  $ (since $\left\vert
I\right\vert \in\mathbb{N}$ and $\left\vert I\right\vert \leq n$). Qed.}.
\end{verlong}

\begin{vershort}
\textbf{(a)} Exercise \ref{exe.prod(ai+bi)} \textbf{(a)} (applied to
$a_{i}=x_{i}$ and $b_{i}=t$) yields%
\begin{align*}
\prod_{i=1}^{n}\left(  x_{i}+t\right)   &  =\underbrace{\sum_{I\subseteq
\left[  n\right]  }}_{\substack{=\sum_{j\in\left\{  0,1,\ldots,n\right\}
}\sum_{\substack{I\subseteq\left[  n\right]  ;\\\left\vert I\right\vert
=j}}\\\text{(since every subset }I\text{ of }\left[  n\right]
\\\text{satisfies }\left\vert I\right\vert \in\left\{  0,1,\ldots,n\right\}
\text{)}}}\left(  \prod_{i\in I}x_{i}\right)  \underbrace{\left(  \prod
_{i\in\left[  n\right]  \setminus I}t\right)  }_{\substack{=t^{\left\vert
\left[  n\right]  \setminus I\right\vert }=t^{\left\vert \left[  n\right]
\right\vert -\left\vert I\right\vert }\\\text{(since }\left\vert \left[
n\right]  \setminus I\right\vert =\left\vert \left[  n\right]  \right\vert
-\left\vert I\right\vert \\\text{(because }I\subseteq\left[  n\right]
\text{))}}}\\
&  =\underbrace{\sum_{j\in\left\{  0,1,\ldots,n\right\}  }}_{=\sum_{j=0}^{n}%
}\ \ \sum_{\substack{I\subseteq\left[  n\right]  ;\\\left\vert I\right\vert
=j}}\left(  \prod_{i\in I}x_{i}\right)  \underbrace{t^{\left\vert \left[
n\right]  \right\vert -\left\vert I\right\vert }}_{\substack{=t^{n-j}%
\\\text{(since }\left\vert \left[  n\right]  \right\vert =n\text{ and
}\left\vert I\right\vert =j\text{)}}}\\
&  =\sum_{j=0}^{n}\ \ \sum_{\substack{I\subseteq\left[  n\right]
;\\\left\vert I\right\vert =j}}\left(  \prod_{i\in I}x_{i}\right)  t^{n-j}.
\end{align*}
Comparing this with%
\begin{align*}
\sum_{j=0}^{n}e_{n-j}\left(  x_{1},x_{2},\ldots,x_{n}\right)  t^{j}  &
=\sum_{j=0}^{n}\underbrace{e_{n-\left(  n-j\right)  }\left(  x_{1}%
,x_{2},\ldots,x_{n}\right)  }_{\substack{=e_{j}\left(  x_{1},x_{2}%
,\ldots,x_{n}\right)  \\=\sum_{\substack{I\subseteq\left[  n\right]
;\\\left\vert I\right\vert =j}}\prod_{i\in I}x_{i}}}t^{n-j}\\
&  \ \ \ \ \ \ \ \ \ \ \left(  \text{here, we have substituted }n-j\text{ for
}j\text{ in the sum}\right) \\
&  =\sum_{j=0}^{n}\ \ \sum_{\substack{I\subseteq\left[  n\right]
;\\\left\vert I\right\vert =j}}\left(  \prod_{i\in I}x_{i}\right)  t^{n-j},
\end{align*}
we obtain%
\[
\prod_{i=1}^{n}\left(  x_{i}+t\right)  =\sum_{j=0}^{n}e_{n-j}\left(
x_{1},x_{2},\ldots,x_{n}\right)  t^{j}.
\]
This proves Proposition \ref{prop.sol.vander-hook.viete} \textbf{(a)}.
\end{vershort}

\begin{verlong}
\textbf{(a)} Exercise \ref{exe.prod(ai+bi)} \textbf{(a)} (applied to
$a_{i}=x_{i}$ and $b_{i}=t$) yields%
\begin{align*}
\prod_{i=1}^{n}\left(  x_{i}+t\right)   &  =\underbrace{\sum_{I\subseteq
\left[  n\right]  }}_{\substack{=\sum_{j\in\left\{  0,1,\ldots,n\right\}
}\sum_{\substack{I\subseteq\left[  n\right]  ;\\\left\vert I\right\vert
=j}}\\\text{(since every subset }I\text{ of }\left[  n\right]
\\\text{satisfies }\left\vert I\right\vert \in\left\{  0,1,\ldots,n\right\}
\text{)}}}\left(  \prod_{i\in I}x_{i}\right)  \underbrace{\left(  \prod
_{i\in\left[  n\right]  \setminus I}t\right)  }_{\substack{=t^{\left\vert
\left[  n\right]  \setminus I\right\vert }=t^{\left\vert \left[  n\right]
\right\vert -\left\vert I\right\vert }\\\text{(since }\left\vert \left[
n\right]  \setminus I\right\vert =\left\vert \left[  n\right]  \right\vert
-\left\vert I\right\vert \\\text{(because }I\subseteq\left[  n\right]
\text{))}}}\\
&  =\underbrace{\sum_{j\in\left\{  0,1,\ldots,n\right\}  }}_{=\sum_{j=0}^{n}%
}\ \ \sum_{\substack{I\subseteq\left[  n\right]  ;\\\left\vert I\right\vert
=j}}\left(  \prod_{i\in I}x_{i}\right)  \underbrace{t^{\left\vert \left[
n\right]  \right\vert -\left\vert I\right\vert }}_{\substack{=t^{n-j}%
\\\text{(since }\left\vert \left[  n\right]  \right\vert =n\text{ and
}\left\vert I\right\vert =j\text{)}}}\\
&  =\sum_{j=0}^{n}\ \ \sum_{\substack{I\subseteq\left[  n\right]
;\\\left\vert I\right\vert =j}}\left(  \prod_{i\in I}x_{i}\right)  t^{n-j}.
\end{align*}
Comparing this with%
\begin{align*}
&  \sum_{j=0}^{n}e_{n-j}\left(  x_{1},x_{2},\ldots,x_{n}\right)  t^{j}\\
&  =\sum_{j=0}^{n}\underbrace{e_{n-\left(  n-j\right)  }\left(  x_{1}%
,x_{2},\ldots,x_{n}\right)  }_{\substack{=e_{j}\left(  x_{1},x_{2}%
,\ldots,x_{n}\right)  \\\text{(since }n-\left(  n-j\right)  =j\text{)}%
}}t^{n-j}\\
&  \ \ \ \ \ \ \ \ \ \ \left(  \text{here, we have substituted }n-j\text{ for
}j\text{ in the sum}\right) \\
&  =\sum_{j=0}^{n}\underbrace{e_{j}\left(  x_{1},x_{2},\ldots,x_{n}\right)
}_{\substack{=\sum_{\substack{I\subseteq\left[  n\right]  ;\\\left\vert
I\right\vert =j}}\prod_{i\in I}x_{i}\\\text{(by the definition of }%
e_{j}\left(  x_{1},x_{2},\ldots,x_{n}\right)  \text{)}}}t^{n-j}\\
&  =\sum_{j=0}^{n}\ \ \sum_{\substack{I\subseteq\left[  n\right]
;\\\left\vert I\right\vert =j}}\left(  \prod_{i\in I}x_{i}\right)  t^{n-j},
\end{align*}
we obtain%
\[
\prod_{i=1}^{n}\left(  x_{i}+t\right)  =\sum_{j=0}^{n}e_{n-j}\left(
x_{1},x_{2},\ldots,x_{n}\right)  t^{j}.
\]
This proves Proposition \ref{prop.sol.vander-hook.viete} \textbf{(a)}.
\end{verlong}

\textbf{(b)} If $I$ is a subset of $\left[  n\right]  $, then we have the
following logical equivalence:%
\begin{equation}
\left(  \left\vert I\right\vert =n\right)  \ \Longleftrightarrow\ \left(
I=\left[  n\right]  \right)  . \label{pf.prop.sol.vander-hook.viete.b.equiv}%
\end{equation}

\begin{vershort}
(The proof of (\ref{pf.prop.sol.vander-hook.viete.b.equiv}) is obvious, since
$\left\vert \left[  n\right]  \right\vert =n$.)
\end{vershort}

\begin{verlong}
[\textit{Proof of (\ref{pf.prop.sol.vander-hook.viete.b.equiv}):} Let $I$ be a
subset of $\left[  n\right]  $.

Assume first that $\left\vert I\right\vert =n$. From $I\subseteq\left[
n\right]  $, we obtain $\left\vert \left[  n\right]  \right\vert =\left\vert
I\right\vert +\left\vert \left[  n\right]  \setminus I\right\vert $ (since the
set $\left[  n\right]  $ is finite). Comparing this with $\left\vert \left[
n\right]  \right\vert =n$, we find $\left\vert I\right\vert +\left\vert
\left[  n\right]  \setminus I\right\vert =n$. Since $\left\vert I\right\vert
=n$, this rewrites as $n+\left\vert \left[  n\right]  \setminus I\right\vert
=n$. Subtracting $n$ from both sides of this equality, we obtain $\left\vert
\left[  n\right]  \setminus I\right\vert =0$. Thus, $\left[  n\right]
\setminus I=\varnothing$. In other words, $\left[  n\right]  \subseteq I$.
Combined with $I\subseteq\left[  n\right]  $, this yields $I=\left[  n\right]
$.

Now, let us forget our assumption that $\left\vert I\right\vert =n$. We thus
have shown that $I=\left[  n\right]  $ under the assumption that $\left\vert
I\right\vert =n$. In other words, we have proven the implication%
\begin{equation}
\left(  \left\vert I\right\vert =n\right)  \ \Longrightarrow\ \left(
I=\left[  n\right]  \right)  .
\label{pf.prop.sol.vander-hook.viete.b.equiv.pf.1}%
\end{equation}

On the other hand, if $I=\left[  n\right]  $, then $\left\vert I\right\vert
=n$ (because if $I=\left[  n\right]  $, then $\left\vert \underbrace{I}%
_{=\left[  n\right]  }\right\vert =\left\vert \left[  n\right]  \right\vert
=n$). In other words, the implication%
\[
\left(  I=\left[  n\right]  \right)  \ \Longrightarrow\ \left(  \left\vert
I\right\vert =n\right)
\]
holds. Combining this implication with
(\ref{pf.prop.sol.vander-hook.viete.b.equiv.pf.1}), we obtain the equivalence
$\left(  \left\vert I\right\vert =n\right)  \ \Longleftrightarrow\ \left(
I=\left[  n\right]  \right)  $. Thus,
(\ref{pf.prop.sol.vander-hook.viete.b.equiv}) is proven.]
\end{verlong}

\begin{vershort}
Now, the definition of $e_{n}\left(  x_{1},x_{2},\ldots,x_{n}\right)  $ yields%
\begin{align*}
e_{n}\left(  x_{1},x_{2},\ldots,x_{n}\right)   &  =\underbrace{\sum
_{\substack{I\subseteq\left[  n\right]  ;\\\left\vert I\right\vert =n}%
}}_{\substack{=\sum_{\substack{I\subseteq\left[  n\right]  ;\\I=\left[
n\right]  }}\\\text{(by the equivalence
(\ref{pf.prop.sol.vander-hook.viete.b.equiv}))}}}\prod_{i\in I}x_{i}%
=\sum_{\substack{I\subseteq\left[  n\right]  ;\\I=\left[  n\right]
}}\ \ \prod_{i\in I}x_{i}\\
&  =\underbrace{\prod_{i\in\left[  n\right]  }}_{=\prod_{i=1}^{n}}%
x_{i}\ \ \ \ \ \ \ \ \ \ \left(  \text{since }\left[  n\right]  \text{ is a
subset of }\left[  n\right]  \right) \\
&  =\prod_{i=1}^{n}x_{i}.
\end{align*}
This proves Proposition \ref{prop.sol.vander-hook.viete} \textbf{(b)}.
\end{vershort}

\begin{verlong}
Now, the definition of $e_{n}\left(  x_{1},x_{2},\ldots,x_{n}\right)  $ yields%
\[
e_{n}\left(  x_{1},x_{2},\ldots,x_{n}\right)  =\underbrace{\sum
_{\substack{I\subseteq\left[  n\right]  ;\\\left\vert I\right\vert =n}%
}}_{\substack{=\sum_{\substack{I\subseteq\left[  n\right]  ;\\I=\left[
n\right]  }}\\\text{(because for every subset }I\text{ of }\left[  n\right]
\text{,}\\\text{the condition }\left(  \left\vert I\right\vert =n\right)
\text{ is equivalent}\\\text{to the condition }\left(  I=\left[  n\right]
\right)  \\\text{(by (\ref{pf.prop.sol.vander-hook.viete.b.equiv})))}}%
}\prod_{i\in I}x_{i}=\sum_{\substack{I\subseteq\left[  n\right]  ;\\I=\left[
n\right]  }}\ \ \prod_{i\in I}x_{i}.
\]
But $\left[  n\right]  $ is a subset of $\left[  n\right]  $. Hence, the sum
$\sum_{\substack{I\subseteq\left[  n\right]  ;\\I=\left[  n\right]  }%
}\prod_{i\in I}x_{i}$ has exactly one addend: namely, the addend for
$I=\left[  n\right]  $. Therefore, this sum simplifies as follows:%
\[
\sum_{\substack{I\subseteq\left[  n\right]  ;\\I=\left[  n\right]  }%
}\ \ \prod_{i\in I}x_{i}=\underbrace{\prod_{i\in\left[  n\right]  }%
}_{\substack{=\prod_{i\in\left\{  1,2,\ldots,n\right\}  }\\\text{(since
}\left[  n\right]  =\left\{  1,2,\ldots,n\right\}  \text{)}}}x_{i}%
=\underbrace{\prod_{i\in\left\{  1,2,\ldots,n\right\}  }}_{=\prod_{i=1}^{n}%
}x_{i}=\prod_{i=1}^{n}x_{i}.
\]
Thus,%
\[
e_{n}\left(  x_{1},x_{2},\ldots,x_{n}\right)  =\sum_{\substack{I\subseteq
\left[  n\right]  ;\\I=\left[  n\right]  }}\ \ \prod_{i\in I}x_{i}=\prod
_{i=1}^{n}x_{i}.
\]
This proves Proposition \ref{prop.sol.vander-hook.viete} \textbf{(b)}.
\end{verlong}

\textbf{(c)} We have $0\in\left\{  0,1,\ldots,n\right\}  $ (since
$n\in\mathbb{N}$). Proposition \ref{prop.sol.vander-hook.viete} \textbf{(a)}
shows that%
\begin{align*}
\prod_{i=1}^{n}\left(  x_{i}+t\right)   &  =\sum_{j=0}^{n}e_{n-j}\left(
x_{1},x_{2},\ldots,x_{n}\right)  t^{j}\\
&  =\underbrace{e_{n-0}\left(  x_{1},x_{2},\ldots,x_{n}\right)  }%
_{\substack{=e_{n}\left(  x_{1},x_{2},\ldots,x_{n}\right)  \\=\prod_{i=1}%
^{n}x_{i}\\\text{(by Proposition \ref{prop.sol.vander-hook.viete}
\textbf{(b)})}}}\underbrace{t^{0}}_{=1}+\underbrace{\sum_{j=1}^{n}%
e_{n-j}\left(  x_{1},x_{2},\ldots,x_{n}\right)  t^{j}}_{\substack{=\sum
_{j=0}^{n-1}e_{n-\left(  j+1\right)  }\left(  x_{1},x_{2},\ldots,x_{n}\right)
t^{j+1}\\\text{(here, we have substituted }j+1\\\text{for }j\text{ in the
sum)}}}\\
&  \ \ \ \ \ \ \ \ \ \ \left(
\begin{array}
[c]{c}%
\text{here, we have split off the addend for }j=0\text{ from the sum}\\
\text{(since }0\in\left\{  0,1,\ldots,n\right\}  \text{)}%
\end{array}
\right) \\
&  =\prod_{i=1}^{n}x_{i}+\sum_{j=0}^{n-1}\underbrace{e_{n-\left(  j+1\right)
}\left(  x_{1},x_{2},\ldots,x_{n}\right)  }_{\substack{=e_{n-1-j}\left(
x_{1},x_{2},\ldots,x_{n}\right)  \\\text{(since }n-\left(  j+1\right)
=n-1-j\text{)}}}\underbrace{t^{j+1}}_{=tt^{j}}\\
&  =\prod_{i=1}^{n}x_{i}+\sum_{j=0}^{n-1}e_{n-1-j}\left(  x_{1},x_{2}%
,\ldots,x_{n}\right)  tt^{j}.
\end{align*}
Solving this equality for $\sum_{j=0}^{n-1}e_{n-1-j}\left(  x_{1},x_{2}%
,\ldots,x_{n}\right)  tt^{j}$, we obtain%
\begin{equation}
\sum_{j=0}^{n-1}e_{n-1-j}\left(  x_{1},x_{2},\ldots,x_{n}\right)  tt^{j}%
=\prod_{i=1}^{n}\left(  x_{i}+t\right)  -\prod_{i=1}^{n}x_{i}.
\label{pf.prop.sol.vander-hook.viete.c.1}%
\end{equation}

But
\begin{align*}
&  t\cdot\underbrace{z_{t}\left(  x_{1},x_{2},\ldots,x_{n}\right)
}_{\substack{=\sum_{j=0}^{n-1}e_{n-1-j}\left(  x_{1},x_{2},\ldots
,x_{n}\right)  t^{j}\\\text{(by the definition of }z_{t}\left(  x_{1}%
,x_{2},\ldots,x_{n}\right)  \text{)}}}\\
&  =t\cdot\sum_{j=0}^{n-1}e_{n-1-j}\left(  x_{1},x_{2},\ldots,x_{n}\right)
t^{j}=\sum_{j=0}^{n-1}e_{n-1-j}\left(  x_{1},x_{2},\ldots,x_{n}\right)
tt^{j}\\
&  =\prod_{i=1}^{n}\left(  x_{i}+t\right)  -\prod_{i=1}^{n}x_{i}%
\ \ \ \ \ \ \ \ \ \ \left(  \text{by (\ref{pf.prop.sol.vander-hook.viete.c.1}%
)}\right)  .
\end{align*}
This proves Proposition \ref{prop.sol.vander-hook.viete} \textbf{(c)}.

\textbf{(d)} Proposition \ref{prop.sol.vander-hook.viete} \textbf{(c)} yields%
\[
t\cdot z_{t}\left(  x_{1},x_{2},\ldots,x_{n}\right)  =\prod_{i=1}^{n}\left(
x_{i}+t\right)  -\prod_{i=1}^{n}x_{i}.
\]
We can divide both sides of this equality by $t$ (since $t$ is invertible). As
a result, we obtain%
\[
z_{t}\left(  x_{1},x_{2},\ldots,x_{n}\right)  =\dfrac{\prod_{i=1}^{n}\left(
x_{i}+t\right)  -\prod_{i=1}^{n}x_{i}}{t}.
\]
This proves Proposition \ref{prop.sol.vander-hook.viete} \textbf{(d)}.
\end{proof}

Our next proposition is crucial in getting a grip on the elements $y_{k}$ in
Proposition \ref{prop.sol.vander-hook.variant-1}:

\begin{proposition}
\label{prop.sol.vander-hook.viete-y}Let $n\in\mathbb{N}$. Let $x_{1}%
,x_{2},\ldots,x_{n}$ be $n$ elements of $\mathbb{K}$. For each $i\in\left\{
1,2,\ldots,n\right\}  $, set $y_{i}=\prod_{\substack{j\in\left\{
1,2,\ldots,n\right\}  ;\\j\neq i}}x_{j}$.

Let $k\in\left\{  1,2,\ldots,n\right\}  $. Then,%
\[
y_{k}=\sum_{j=0}^{n-1}\left(  -1\right)  ^{n-1-j}e_{j}\left(  x_{1}%
,x_{2},\ldots,x_{n}\right)  x_{k}^{n-1-j}.
\]

\end{proposition}

\begin{vershort}
Our proof of Proposition \ref{prop.sol.vander-hook.viete-y} will rely on a
basic fact about sets:

\begin{proposition}
\label{prop.sol.vander-hook.viete-y.short.Pm.lem}For every set $T$ and every
$j\in\mathbb{N}$, we let $\mathcal{P}_{j}\left(  T\right)  $ denote the set of
all $j$-element subsets of $T$.

Let $S$ be a set. Let $s\in S$. Let $m$ be a positive integer. Then:

\begin{itemize}
\item We have $\mathcal{P}_{m}\left(  S\setminus\left\{  s\right\}  \right)
\subseteq\mathcal{P}_{m}\left(  S\right)  $.

\item The map%
\begin{align*}
\mathcal{P}_{m-1}\left(  S\setminus\left\{  s\right\}  \right)   &
\rightarrow\mathcal{P}_{m}\left(  S\right)  \setminus\mathcal{P}_{m}\left(
S\setminus\left\{  s\right\}  \right)  ,\\
U  &  \mapsto U\cup\left\{  s\right\}
\end{align*}
is well-defined and a bijection.
\end{itemize}
\end{proposition}

Roughly speaking, Proposition \ref{prop.sol.vander-hook.viete-y.short.Pm.lem}
claims that if $s$ is an element of a set $S$, and if $m$ is a positive
integer, then:

\begin{itemize}
\item all $m$-element subsets of $S\setminus\left\{  s\right\}  $ are
$m$-element subsets of $S$ as well;

\item the $m$-element subsets of $S$ which are \textbf{not} $m$-element
subsets of $S\setminus\left\{  s\right\}  $ are in bijection with the $\left(
m-1\right)  $-element subsets of $S\setminus\left\{  s\right\}  $; this
bijection sends an $\left(  m-1\right)  $-element subset $U$ of $S\setminus
\left\{  s\right\}  $ to the $m$-element subset $U\cup\left\{  s\right\}  $ of
$S$.
\end{itemize}

Restated this way, Proposition \ref{prop.sol.vander-hook.viete-y.short.Pm.lem}
should be intuitively clear. A rigorous proof is not hard to
give\footnote{Here is an outline: First, show that the two maps%
\begin{align*}
\mathcal{P}_{m-1}\left(  S\setminus\left\{  s\right\}  \right)   &
\rightarrow\mathcal{P}_{m}\left(  S\right)  \setminus\mathcal{P}_{m}\left(
S\setminus\left\{  s\right\}  \right)  ,\\
U  &  \mapsto U\cup\left\{  s\right\}
\end{align*}
and%
\begin{align*}
\mathcal{P}_{m}\left(  S\right)  \setminus\mathcal{P}_{m}\left(
S\setminus\left\{  s\right\}  \right)   &  \rightarrow\mathcal{P}_{m-1}\left(
S\setminus\left\{  s\right\}  \right)  ,\\
U  &  \mapsto U\setminus\left\{  s\right\}
\end{align*}
are well-defined. Then, show that these maps are mutually inverse. Conclude
that the first of these two maps is invertible, i.e., is a bijection.}.

\begin{proof}
[Proof of Proposition \ref{prop.sol.vander-hook.viete-y}.]We shall use the
notations $\left[  n\right]  $, $\mathcal{P}\left(  \left[  n\right]  \right)
$ and $\sum_{\substack{I\subseteq\left[  n\right]  ;\\\left\vert I\right\vert
=j}}$ as in Definition \ref{def.sol.vander-hook.elsyms}. We shall furthermore
use the notation $\mathcal{P}_{j}\left(  T\right)  $ from Proposition
\ref{prop.sol.vander-hook.viete-y.short.Pm.lem}.

Let $K=\left[  n\right]  \setminus\left\{  k\right\}  $. Note that
$k\in\left\{  1,2,\ldots,n\right\}  =\left[  n\right]  $, so that $\left\vert
\left[  n\right]  \setminus\left\{  k\right\}  \right\vert
=\underbrace{\left\vert \left[  n\right]  \right\vert }_{=n}-1=n-1$. From
$K=\left[  n\right]  \setminus\left\{  k\right\}  $, we obtain $\left\vert
K\right\vert =\left\vert \left[  n\right]  \setminus\left\{  k\right\}
\right\vert =n-1$.

For every $j\in\mathbb{N}$, define an element $f_{j}\in\mathbb{K}$ by%
\[
f_{j}=\sum_{\substack{I\subseteq K;\\\left\vert I\right\vert =j}%
}\ \ \prod_{i\in I}x_{i}.
\]
Then,%
\begin{equation}
f_{0}=1 \label{pf.prop.sol.vander-hook.viete-y.short.f0=}%
\end{equation}
\footnote{\textit{Proof of (\ref{pf.prop.sol.vander-hook.viete-y.short.f0=}):}
There exists exactly one subset $I$ of $K$ satisfying $\left\vert I\right\vert
=0$: namely, the subset $\varnothing$. Hence, the sum $\sum
_{\substack{I\subseteq K;\\\left\vert I\right\vert =0}}\ \ \prod_{i\in I}%
x_{i}$ has only one addend: namely, the addend for $I=\varnothing$. Thus, this
sum simplifies to $\sum_{\substack{I\subseteq K;\\\left\vert I\right\vert
=0}}\ \ \prod_{i\in I}x_{i}=\prod_{i\in\varnothing}x_{i}=\left(  \text{empty
product}\right)  =1$.
\par
Now, the definition of $f_{0}$ yields $f_{0}=\sum_{\substack{I\subseteq
K;\\\left\vert I\right\vert =0}}\ \ \prod_{i\in I}x_{i}=1$. This proves
(\ref{pf.prop.sol.vander-hook.viete-y.short.f0=}).}.

Next, we shall show that%
\begin{equation}
f_{n-1}=y_{k}. \label{pf.prop.sol.vander-hook.viete-y.short.fn-1=}%
\end{equation}

[\textit{Proof of (\ref{pf.prop.sol.vander-hook.viete-y.short.fn-1=}):} We
have $\left\vert K\right\vert =n-1$. Hence, there exists only one $\left(
n-1\right)  $-element subset $I$ of $K$: namely, the set $K$ itself. In other
words, there exists only one subset $I$ of $K$ satisfying $\left\vert
I\right\vert =n-1$: namely, the set $K$ itself. Hence, the sum $\sum
_{\substack{I\subseteq K;\\\left\vert I\right\vert =n-1}}\prod_{i\in I}x_{i}$
has only one addend: namely, the addend for $I=K$. Thus, this sum simplifies
to%
\begin{align*}
\sum_{\substack{I\subseteq K;\\\left\vert I\right\vert =n-1}}\prod_{i\in
I}x_{i}  &  =\prod_{i\in K}x_{i}=\prod_{j\in K}x_{j}=\underbrace{\prod
_{j\in\left\{  1,2,\ldots,n\right\}  \setminus\left\{  k\right\}  }}%
_{=\prod_{\substack{j\in\left\{  1,2,\ldots,n\right\}  ;\\j\neq k}}}x_{j}\\
&  \ \ \ \ \ \ \ \ \ \ \left(  \text{since }K=\underbrace{\left[  n\right]
}_{=\left\{  1,2,\ldots,n\right\}  }\setminus\left\{  k\right\}  =\left\{
1,2,\ldots,n\right\}  \setminus\left\{  k\right\}  \right) \\
&  =\prod_{\substack{j\in\left\{  1,2,\ldots,n\right\}  ;\\j\neq k}}x_{j}.
\end{align*}
Now, the definition of $f_{n-1}$ yields%
\[
f_{n-1}=\sum_{\substack{I\subseteq K;\\\left\vert I\right\vert =n-1}%
}\prod_{i\in I}x_{i}=\prod_{\substack{j\in\left\{  1,2,\ldots,n\right\}
;\\j\neq k}}x_{j}.
\]
Comparing this with%
\[
y_{k}=\prod_{\substack{j\in\left\{  1,2,\ldots,n\right\}  ;\\j\neq k}%
}x_{j}\ \ \ \ \ \ \ \ \ \ \left(  \text{by the definition of }y_{k}\right)  ,
\]
we find $f_{n-1}=y_{k}$. Thus,
(\ref{pf.prop.sol.vander-hook.viete-y.short.fn-1=}) is proven.]

But%
\begin{equation}
e_{0}\left(  x_{1},x_{2},\ldots,x_{n}\right)  =1
\label{pf.prop.sol.vander-hook.viete-y.short.e0=}%
\end{equation}
\footnote{\textit{Proof of (\ref{pf.prop.sol.vander-hook.viete-y.short.e0=}):}
There exists exactly one subset $I$ of $\left[  n\right]  $ satisfying
$\left\vert I\right\vert =0$: namely, the subset $\varnothing$. Hence, the sum
$\sum_{\substack{I\subseteq\left[  n\right]  ;\\\left\vert I\right\vert
=0}}\ \ \prod_{i\in I}x_{i}$ has only one addend: namely, the addend for
$I=\varnothing$. Thus, this sum simplifies to $\sum_{\substack{I\subseteq
\left[  n\right]  ;\\\left\vert I\right\vert =0}}\ \ \prod_{i\in I}x_{i}%
=\prod_{i\in\varnothing}x_{i}=\left(  \text{empty product}\right)  =1$.
\par
Now, the definition of $e_{0}\left(  x_{1},x_{2},\ldots,x_{n}\right)  $ yields
$e_{0}\left(  x_{1},x_{2},\ldots,x_{n}\right)  =\sum_{\substack{I\subseteq
\left[  n\right]  ;\\\left\vert I\right\vert =0}}\ \ \prod_{i\in I}x_{i}=1$.
This proves (\ref{pf.prop.sol.vander-hook.viete-y.short.e0=}).}.

We shall now show that%
\begin{equation}
e_{m}\left(  x_{1},x_{2},\ldots,x_{n}\right)  =f_{m}+x_{k}f_{m-1}
\label{pf.prop.sol.vander-hook.viete-y.short.main}%
\end{equation}
for every positive integer $m$.

[\textit{Proof of (\ref{pf.prop.sol.vander-hook.viete-y.short.main}):} Let $m$
be a positive integer. The definition of $f_{m}$ yields%
\begin{equation}
f_{m}=\underbrace{\sum_{\substack{I\subseteq K;\\\left\vert I\right\vert =m}%
}}_{\substack{=\sum_{\substack{I\text{ is an }m\text{-element}\\\text{subset
of }K}}=\sum_{I\in\mathcal{P}_{m}\left(  K\right)  }\\\text{(since the set of
all }m\text{-element}\\\text{subsets of }K\text{ is }\mathcal{P}_{m}\left(
K\right)  \text{)}}}\prod_{i\in I}x_{i}=\sum_{I\in\mathcal{P}_{m}\left(
K\right)  }\ \ \prod_{i\in I}x_{i}.
\label{pf.prop.sol.vander-hook.viete-y.short.main.pf.fm=}%
\end{equation}
The same argument (but applied to $m-1$ instead of $m$) yields%
\begin{equation}
f_{m-1}=\sum_{I\in\mathcal{P}_{m-1}\left(  K\right)  }\ \ \prod_{i\in I}x_{i}.
\label{pf.prop.sol.vander-hook.viete-y.short.main.pf.fm-1=}%
\end{equation}

Recall that $k\in\left[  n\right]  $. Hence, we can apply Proposition
\ref{prop.sol.vander-hook.viete-y.short.Pm.lem} to $\left[  n\right]  $ and
$k$ instead of $S$ and $s$. As a result, we obtain the following two results:

\begin{itemize}
\item We have $\mathcal{P}_{m}\left(  \left[  n\right]  \setminus\left\{
k\right\}  \right)  \subseteq\mathcal{P}_{m}\left(  \left[  n\right]  \right)
$.

\item The map%
\begin{align*}
\mathcal{P}_{m-1}\left(  \left[  n\right]  \setminus\left\{  k\right\}
\right)   &  \rightarrow\mathcal{P}_{m}\left(  \left[  n\right]  \right)
\setminus\mathcal{P}_{m}\left(  \left[  n\right]  \setminus\left\{  k\right\}
\right)  ,\\
U  &  \mapsto U\cup\left\{  k\right\}
\end{align*}
is well-defined and a bijection.
\end{itemize}

Since $\left[  n\right]  \setminus\left\{  k\right\}  =K$, these two results
can be rewritten as follows:

\begin{itemize}
\item We have $\mathcal{P}_{m}\left(  K\right)  \subseteq\mathcal{P}%
_{m}\left(  \left[  n\right]  \right)  $.

\item The map%
\begin{align*}
\mathcal{P}_{m-1}\left(  K\right)   &  \rightarrow\mathcal{P}_{m}\left(
\left[  n\right]  \right)  \setminus\mathcal{P}_{m}\left(  K\right)  ,\\
U  &  \mapsto U\cup\left\{  k\right\}
\end{align*}
is well-defined and a bijection.
\end{itemize}

Furthermore, every $I\in\mathcal{P}_{m-1}\left(  K\right)  $ satisfies%
\begin{equation}
\prod_{i\in I\cup\left\{  k\right\}  }x_{i}=x_{k}\prod_{i\in I}x_{i}
\label{pf.prop.sol.vander-hook.viete-y.short.main.pf.2}%
\end{equation}
\footnote{\textit{Proof of
(\ref{pf.prop.sol.vander-hook.viete-y.short.main.pf.2}):} Let $I\in
\mathcal{P}_{m-1}\left(  K\right)  $. Thus, $I$ is an $\left(  m-1\right)
$-element subset of $K$ (since $\mathcal{P}_{m-1}\left(  K\right)  $ is
defined as the set of all $\left(  m-1\right)  $-element subsets of $K$). In
particular, $I$ is a subset of $K$. Thus, $I\subseteq K=\left[  n\right]
\setminus\left\{  k\right\}  $. Thus, $k\notin I$. Hence, the sets $I$ and
$\left\{  k\right\}  $ are disjoint. Thus,%
\[
\prod_{i\in I\cup\left\{  k\right\}  }x_{i}=\left(  \prod_{i\in I}%
x_{i}\right)  \underbrace{\left(  \prod_{i\in\left\{  k\right\}  }%
x_{i}\right)  }_{=x_{k}}=\left(  \prod_{i\in I}x_{i}\right)  x_{k}=x_{k}%
\prod_{i\in I}x_{i}.
\]
This proves (\ref{pf.prop.sol.vander-hook.viete-y.short.main.pf.2}).}. Now,%
\begin{align}
&  \sum_{I\in\mathcal{P}_{m}\left(  \left[  n\right]  \right)  \setminus
\mathcal{P}_{m}\left(  K\right)  }\ \ \prod_{i\in I}x_{i}\nonumber\\
&  =\sum_{U\in\mathcal{P}_{m-1}\left(  K\right)  }\ \ \prod_{i\in
U\cup\left\{  k\right\}  }x_{i}\nonumber\\
&  \ \ \ \ \ \ \ \ \ \ \left(
\begin{array}
[c]{c}%
\text{here, we have substituted }U\cup\left\{  k\right\}  \text{ for }I\text{
in the sum, since}\\
\text{the map }\mathcal{P}_{m-1}\left(  K\right)  \rightarrow\mathcal{P}%
_{m}\left(  \left[  n\right]  \right)  \setminus\mathcal{P}_{m}\left(
K\right)  ,\ U\mapsto U\cup\left\{  k\right\} \\
\text{is a bijection}%
\end{array}
\right) \nonumber\\
&  =\sum_{I\in\mathcal{P}_{m-1}\left(  K\right)  }\ \ \underbrace{\prod_{i\in
I\cup\left\{  k\right\}  }x_{i}}_{\substack{=x_{k}\prod_{i\in I}%
x_{i}\\\text{(by (\ref{pf.prop.sol.vander-hook.viete-y.short.main.pf.2}))}%
}}\ \ \ \ \ \ \ \ \ \ \left(
\begin{array}
[c]{c}%
\text{here, we have renamed the}\\
\text{summation index }U\text{ as }I
\end{array}
\right) \nonumber\\
&  =\sum_{I\in\mathcal{P}_{m-1}\left(  K\right)  }x_{k}\prod_{i\in I}%
x_{i}=x_{k}\sum_{I\in\mathcal{P}_{m-1}\left(  K\right)  }\ \ \prod_{i\in
I}x_{i}. \label{pf.prop.sol.vander-hook.viete-y.short.main.pf.2b}%
\end{align}

But the definition of $e_{m}\left(  x_{1},x_{2},\ldots,x_{n}\right)  $ yields
\begin{align*}
e_{m}\left(  x_{1},x_{2},\ldots,x_{n}\right)   &  =\underbrace{\sum
_{\substack{I\subseteq\left[  n\right]  ;\\\left\vert I\right\vert =m}%
}}_{\substack{=\sum_{\substack{I\text{ is an }m\text{-element}\\\text{subset
of }\left[  n\right]  }}=\sum_{I\in\mathcal{P}_{m}\left(  \left[  n\right]
\right)  }\\\text{(since the set of all }m\text{-element}\\\text{subsets of
}\left[  n\right]  \text{ is }\mathcal{P}_{m}\left(  \left[  n\right]
\right)  \text{)}}}\prod_{i\in I}x_{i}=\sum_{I\in\mathcal{P}_{m}\left(
\left[  n\right]  \right)  }\ \ \prod_{i\in I}x_{i}\\
&  =\underbrace{\sum_{\substack{I\in\mathcal{P}_{m}\left(  \left[  n\right]
\right)  ;\\I\in\mathcal{P}_{m}\left(  K\right)  }}}_{\substack{=\sum
_{I\in\mathcal{P}_{m}\left(  K\right)  }\\\text{(since }\mathcal{P}_{m}\left(
K\right)  \subseteq\mathcal{P}_{m}\left(  \left[  n\right]  \right)  \text{)}%
}}\prod_{i\in I}x_{i}+\underbrace{\sum_{\substack{I\in\mathcal{P}_{m}\left(
\left[  n\right]  \right)  ;\\I\notin\mathcal{P}_{m}\left(  K\right)  }%
}}_{=\sum_{I\in\mathcal{P}_{m}\left(  \left[  n\right]  \right)
\setminus\mathcal{P}_{m}\left(  K\right)  }}\prod_{i\in I}x_{i}\\
&  \ \ \ \ \ \ \ \ \ \ \left(
\begin{array}
[c]{c}%
\text{since each }I\in\mathcal{P}_{m}\left(  \left[  n\right]  \right)  \text{
satisfies either }I\in\mathcal{P}_{m}\left(  K\right) \\
\text{or }I\notin\mathcal{P}_{m}\left(  K\right)  \text{ (but not both)}%
\end{array}
\right) \\
&  =\sum_{I\in\mathcal{P}_{m}\left(  K\right)  }\ \ \prod_{i\in I}%
x_{i}+\underbrace{\sum_{I\in\mathcal{P}_{m}\left(  \left[  n\right]  \right)
\setminus\mathcal{P}_{m}\left(  K\right)  }\ \ \prod_{i\in I}x_{i}%
}_{\substack{=x_{k}\sum_{I\in\mathcal{P}_{m-1}\left(  K\right)  }%
\ \ \prod_{i\in I}x_{i}\\\text{(by
(\ref{pf.prop.sol.vander-hook.viete-y.short.main.pf.2b}))}}}\\
&  =\underbrace{\sum_{I\in\mathcal{P}_{m}\left(  K\right)  }\ \ \prod_{i\in
I}x_{i}}_{\substack{=f_{m}\\\text{(by
(\ref{pf.prop.sol.vander-hook.viete-y.short.main.pf.fm=}))}}}+x_{k}%
\underbrace{\sum_{I\in\mathcal{P}_{m-1}\left(  K\right)  }\ \ \prod_{i\in
I}x_{i}}_{\substack{=f_{m-1}\\\text{(by
(\ref{pf.prop.sol.vander-hook.viete-y.short.main.pf.fm-1=}))}}}\\
&  =f_{m}+x_{k}f_{m-1}.
\end{align*}
This proves (\ref{pf.prop.sol.vander-hook.viete-y.short.main}).]

Recall that $k\in\left\{  1,2,\ldots,n\right\}  $. Hence, $1\leq k\leq n$, so
that $n\geq1$, so that $n-1\in\mathbb{N}$. Hence, $0\in\left\{  0,1,\ldots
,n-1\right\}  $ and $n-1\in\left\{  0,1,\ldots,n-1\right\}  $.

Now,%
\begin{align}
&  \sum_{j=0}^{n-1}\left(  -1\right)  ^{n-1-j}e_{j}\left(  x_{1},x_{2}%
,\ldots,x_{n}\right)  x_{k}^{n-1-j}\nonumber\\
&  =\underbrace{\left(  -1\right)  ^{n-1-0}}_{=\left(  -1\right)  ^{n-1}%
}\underbrace{e_{0}\left(  x_{1},x_{2},\ldots,x_{n}\right)  }%
_{\substack{=1\\\text{(by (\ref{pf.prop.sol.vander-hook.viete-y.short.e0=}))}%
}}\underbrace{x_{k}^{n-1-0}}_{=x_{k}^{n-1}}+\sum_{j=1}^{n-1}\left(  -1\right)
^{n-1-j}\underbrace{e_{j}\left(  x_{1},x_{2},\ldots,x_{n}\right)
}_{\substack{=f_{j}+x_{k}f_{j-1}\\\text{(by
(\ref{pf.prop.sol.vander-hook.viete-y.short.main}) (applied to }m=j\text{))}%
}}x_{k}^{n-1-j}\nonumber\\
&  \ \ \ \ \ \ \ \ \ \ \ \ \ \ \ \ \ \ \ \ \left(
\begin{array}
[c]{c}%
\text{here, we have split off the addend for }j=0\text{ from the sum,}\\
\text{since }0\in\left\{  0,1,\ldots,n-1\right\}
\end{array}
\right) \nonumber\\
&  =\left(  -1\right)  ^{n-1}x_{k}^{n-1}+\underbrace{\sum_{j=1}^{n-1}\left(
-1\right)  ^{n-1-j}\left(  f_{j}+x_{k}f_{j-1}\right)  x_{k}^{n-1-j}}%
_{=\sum_{j=1}^{n-1}\left(  -1\right)  ^{n-1-j}f_{j}x_{k}^{n-1-j}+\sum
_{j=1}^{n-1}\left(  -1\right)  ^{n-1-j}x_{k}f_{j-1}x_{k}^{n-1-j}}\nonumber\\
&  =\left(  -1\right)  ^{n-1}x_{k}^{n-1}+\sum_{j=1}^{n-1}\left(  -1\right)
^{n-1-j}f_{j}x_{k}^{n-1-j}+\sum_{j=1}^{n-1}\left(  -1\right)  ^{n-1-j}%
x_{k}f_{j-1}x_{k}^{n-1-j}. \label{pf.prop.sol.vander-hook.viete-y.short.5}%
\end{align}

But
\begin{align}
&  \sum_{j=1}^{n-1}\left(  -1\right)  ^{n-1-j}\underbrace{x_{k}f_{j-1}%
}_{=f_{j-1}x_{k}}x_{k}^{n-1-j}\nonumber\\
&  =\sum_{j=1}^{n-1}\left(  -1\right)  ^{n-1-j}f_{j-1}\underbrace{x_{k}%
x_{k}^{n-1-j}}_{=x_{k}^{\left(  n-1-j\right)  +1}=x_{k}^{n-j}}=\sum
_{j=1}^{n-1}\left(  -1\right)  ^{n-1-j}f_{j-1}x_{k}^{n-j}\nonumber\\
&  =\sum_{j=0}^{\left(  n-1\right)  -1}\underbrace{\left(  -1\right)
^{n-1-\left(  j+1\right)  }}_{\substack{=\left(  -1\right)  ^{n-j}%
\\\text{(since }n-1-\left(  j+1\right)  =n-j-2\equiv n-j\operatorname{mod}%
2\text{)}}}\underbrace{f_{\left(  j+1\right)  -1}}_{=f_{j}}\underbrace{x_{k}%
^{n-\left(  j+1\right)  }}_{\substack{=x_{k}^{n-1-j}\\\text{(since }n-\left(
j+1\right)  =n-1-j\text{)}}}\nonumber\\
&  \ \ \ \ \ \ \ \ \ \ \ \ \ \ \ \ \ \ \ \ \left(  \text{here, we have
substituted }j+1\text{ for }j\text{ in the sum}\right) \nonumber\\
&  =\sum_{j=0}^{\left(  n-1\right)  -1}\underbrace{\left(  -1\right)  ^{n-j}%
}_{\substack{=-\left(  -1\right)  ^{n-j-1}=-\left(  -1\right)  ^{n-1-j}%
\\\text{(since }n-j-1=n-1-j\text{)}}}f_{j}x_{k}^{n-1-j}=\sum_{j=0}^{\left(
n-1\right)  -1}\left(  -\left(  -1\right)  ^{n-1-j}\right)  f_{j}x_{k}%
^{n-1-j}\nonumber\\
&  =-\sum_{j=0}^{\left(  n-1\right)  -1}\left(  -1\right)  ^{n-1-j}f_{j}%
x_{k}^{n-1-j}. \label{pf.prop.sol.vander-hook.viete-y.short.6a}%
\end{align}
Also,%
\begin{align*}
&  \sum_{j=0}^{n-1}\left(  -1\right)  ^{n-1-j}f_{j}x_{k}^{n-1-j}\\
&  =\underbrace{\left(  -1\right)  ^{n-1-0}}_{=\left(  -1\right)  ^{n-1}%
}\underbrace{f_{0}}_{\substack{=1\\\text{(by
(\ref{pf.prop.sol.vander-hook.viete-y.short.f0=}))}}}\underbrace{x_{k}%
^{n-1-0}}_{=x_{k}^{n-1}}+\sum_{j=1}^{n-1}\left(  -1\right)  ^{n-1-j}f_{j}%
x_{k}^{n-1-j}\\
&  \ \ \ \ \ \ \ \ \ \ \ \ \ \ \ \ \ \ \ \ \left(
\begin{array}
[c]{c}%
\text{here, we have split off the addend for }j=0\text{ from the sum,}\\
\text{since }0\in\left\{  0,1,\ldots,n-1\right\}
\end{array}
\right) \\
&  =\left(  -1\right)  ^{n-1}x_{k}^{n-1}+\sum_{j=1}^{n-1}\left(  -1\right)
^{n-1-j}f_{j}x_{k}^{n-1-j},
\end{align*}
so that%
\begin{align}
&  \left(  -1\right)  ^{n-1}x_{k}^{n-1}+\sum_{j=1}^{n-1}\left(  -1\right)
^{n-1-j}f_{j}x_{k}^{n-1-j}\nonumber\\
&  =\sum_{j=0}^{n-1}\left(  -1\right)  ^{n-1-j}f_{j}x_{k}^{n-1-j}\nonumber\\
&  =\sum_{j=0}^{\left(  n-1\right)  -1}\left(  -1\right)  ^{n-1-j}f_{j}%
x_{k}^{n-1-j}+\underbrace{\left(  -1\right)  ^{n-1-\left(  n-1\right)  }%
}_{=\left(  -1\right)  ^{0}=1}\underbrace{f_{n-1}}_{\substack{=y_{k}%
\\\text{(by (\ref{pf.prop.sol.vander-hook.viete-y.short.fn-1=}))}%
}}\underbrace{x_{k}^{n-1-\left(  n-1\right)  }}_{=x_{k}^{0}=1}\nonumber\\
&  \ \ \ \ \ \ \ \ \ \ \ \ \ \ \ \ \ \ \ \ \left(
\begin{array}
[c]{c}%
\text{here, we have split off the addend for }j=n-1\text{ from the sum,}\\
\text{since }n-1\in\left\{  0,1,\ldots,n-1\right\}
\end{array}
\right) \nonumber\\
&  =\sum_{j=0}^{\left(  n-1\right)  -1}\left(  -1\right)  ^{n-1-j}f_{j}%
x_{k}^{n-1-j}+y_{k}. \label{pf.prop.sol.vander-hook.viete-y.short.6b}%
\end{align}

Now, (\ref{pf.prop.sol.vander-hook.viete-y.short.5}) becomes%
\begin{align*}
&  \sum_{j=0}^{n-1}\left(  -1\right)  ^{n-1-j}e_{j}\left(  x_{1},x_{2}%
,\ldots,x_{n}\right)  x_{k}^{n-1-j}\\
&  =\underbrace{\left(  -1\right)  ^{n-1}x_{k}^{n-1}+\sum_{j=1}^{n-1}\left(
-1\right)  ^{n-1-j}f_{j}x_{k}^{n-1-j}}_{\substack{=\sum_{j=0}^{\left(
n-1\right)  -1}\left(  -1\right)  ^{n-1-j}f_{j}x_{k}^{n-1-j}+y_{k}\\\text{(by
(\ref{pf.prop.sol.vander-hook.viete-y.short.6b}))}}}+\underbrace{\sum
_{j=1}^{n-1}\left(  -1\right)  ^{n-1-j}x_{k}f_{j-1}x_{k}^{n-1-j}%
}_{\substack{=-\sum_{j=0}^{\left(  n-1\right)  -1}\left(  -1\right)
^{n-1-j}f_{j}x_{k}^{n-1-j}\\\text{(by
(\ref{pf.prop.sol.vander-hook.viete-y.short.6a}))}}}\\
&  =\sum_{j=0}^{\left(  n-1\right)  -1}\left(  -1\right)  ^{n-1-j}f_{j}%
x_{k}^{n-1-j}+y_{k}+\left(  -\sum_{j=0}^{\left(  n-1\right)  -1}\left(
-1\right)  ^{n-1-j}f_{j}x_{k}^{n-1-j}\right) \\
&  =y_{k}.
\end{align*}
This proves Proposition \ref{prop.sol.vander-hook.viete-y}.
\end{proof}
\end{vershort}

\begin{verlong}
\begin{proof}
[Proof of Proposition \ref{prop.sol.vander-hook.viete-y}.]We shall use the
notations $\left[  n\right]  $, $\mathcal{P}\left(  \left[  n\right]  \right)
$ and $\sum_{\substack{I\subseteq\left[  n\right]  ;\\\left\vert I\right\vert
=j}}$ as in Definition \ref{def.sol.vander-hook.elsyms}. Let $K=\left[
n\right]  \setminus\left\{  k\right\}  $. Note that $k\in\left\{
1,2,\ldots,n\right\}  =\left[  n\right]  $, so that $\left\vert \left[
n\right]  \setminus\left\{  k\right\}  \right\vert =\underbrace{\left\vert
\left[  n\right]  \right\vert }_{=n}-1=n-1$. From $K=\left[  n\right]
\setminus\left\{  k\right\}  $, we obtain $\left\vert K\right\vert =\left\vert
\left[  n\right]  \setminus\left\{  k\right\}  \right\vert =n-1$. Thus, $K$ is
a finite set.

We shall use the notation introduced in Definition
\ref{def.sol.prop.binom.subsets.Pm}.

For every $j\in\mathbb{N}$, define an element $f_{j}\in\mathbb{K}$ by%
\[
f_{j}=\sum_{I\in\mathcal{P}_{j}\left(  K\right)  }\ \ \prod_{i\in I}x_{i}.
\]
Then,%
\begin{equation}
f_{0}=1 \label{pf.prop.sol.vander-hook.viete-y.f0=}%
\end{equation}
\footnote{\textit{Proof of (\ref{pf.prop.sol.vander-hook.viete-y.f0=}):}
Proposition \ref{prop.sol.prop.binom.subsets.Pm.lem} \textbf{(a)} (applied to
$S=K$) yields $\mathcal{P}_{0}\left(  K\right)  =\left\{  \varnothing\right\}
$. Now, the definition of $f_{0}$ yields%
\begin{align*}
f_{0}  &  =\sum_{I\in\mathcal{P}_{0}\left(  K\right)  }\ \ \prod_{i\in I}%
x_{i}=\sum_{I\in\left\{  \varnothing\right\}  }\ \ \prod_{i\in I}%
x_{i}\ \ \ \ \ \ \ \ \ \ \left(  \text{since }\mathcal{P}_{0}\left(  K\right)
=\left\{  \varnothing\right\}  \right) \\
&  =\prod_{i\in\varnothing}x_{i}=\left(  \text{empty product}\right)  =1.
\end{align*}
This proves (\ref{pf.prop.sol.vander-hook.viete-y.f0=}).}.

Next, we shall show that%
\begin{equation}
f_{n-1}=y_{k}. \label{pf.prop.sol.vander-hook.viete-y.fn-1=}%
\end{equation}

[\textit{Proof of (\ref{pf.prop.sol.vander-hook.viete-y.fn-1=}):} We have
$\left\{  K\right\}  \subseteq\mathcal{P}_{n-1}\left(  K\right)
$\ \ \ \ \footnote{\textit{Proof.} We know that $\mathcal{P}_{n-1}\left(
K\right)  $ is the set of all $\left(  n-1\right)  $-element subsets of $K$
(by the definition of $\mathcal{P}_{n-1}\left(  K\right)  $).
\par
But $K$ is an $\left(  n-1\right)  $-element set (since $\left\vert
K\right\vert =n-1$) and a subset of $K$ (obviously). Hence, $K$ is an $\left(
n-1\right)  $-element subset of $K$. In other words, $K\in\mathcal{P}%
_{n-1}\left(  K\right)  $ (since $\mathcal{P}_{n-1}\left(  K\right)  $ is the
set of all $\left(  n-1\right)  $-element subsets of $K$). Thus, $\left\{
K\right\}  \subseteq\mathcal{P}_{n-1}\left(  K\right)  $. Qed.} and
$\mathcal{P}_{n-1}\left(  K\right)  \subseteq\left\{  K\right\}
$\ \ \ \ \footnote{\textit{Proof.} We know that $\mathcal{P}_{n-1}\left(
K\right)  $ is the set of all $\left(  n-1\right)  $-element subsets of $K$
(by the definition of $\mathcal{P}_{n-1}\left(  K\right)  $).
\par
Let $I\in\mathcal{P}_{n-1}\left(  K\right)  $. Thus, $I$ is an $\left(
n-1\right)  $-element subset of $K$ (since $\mathcal{P}_{n-1}\left(  K\right)
$ is the set of all $\left(  n-1\right)  $-element subsets of $K$). In other
words, $I$ is a subset of $K$ and satisfies $\left\vert I\right\vert =n-1$.
\par
We have $I\subseteq K$ (since $I$ is a subset of $K$) and thus $\left\vert
K\right\vert =\left\vert I\right\vert +\left\vert K\setminus I\right\vert $
(since $K$ is a finite set). Compared with $\left\vert K\right\vert =n-1$,
this yields $n-1=\underbrace{\left\vert I\right\vert }_{=n-1}+\left\vert
K\setminus I\right\vert =n-1+\left\vert K\setminus I\right\vert $. Subtracting
$n-1$ from both sides of this equality, we find $0=\left\vert K\setminus
I\right\vert $. Thus, $\left\vert K\setminus I\right\vert =0$, so that
$K\setminus I=\varnothing$. In other words, $K\subseteq I$. Combined with
$I\subseteq K$, this yields $I=K$. Thus, $I\in\left\{  K\right\}  $.
\par
Now, forget that we fixed $I$. We thus have shown that every $I\in
\mathcal{P}_{n-1}\left(  K\right)  $ satisfies $I\in\left\{  K\right\}  $. In
other words, $\mathcal{P}_{n-1}\left(  K\right)  \subseteq\left\{  K\right\}
$. Qed.}. Combining these two relations, we obtain $\mathcal{P}_{n-1}\left(
K\right)  =\left\{  K\right\}  $. Now, the definition of $f_{n-1}$ yields%
\begin{align*}
f_{n-1}  &  =\sum_{I\in\mathcal{P}_{n-1}\left(  K\right)  }\ \ \prod_{i\in
I}x_{i}=\sum_{I\in\left\{  K\right\}  }\ \ \prod_{i\in I}x_{i}%
\ \ \ \ \ \ \ \ \ \ \left(  \text{since }\mathcal{P}_{n-1}\left(  K\right)
=\left\{  K\right\}  \right) \\
&  =\underbrace{\prod_{i\in K}}_{\substack{=\prod_{i\in\left[  n\right]
\setminus\left\{  k\right\}  }\\\text{(since }K=\left[  n\right]
\setminus\left\{  k\right\}  \text{)}}}x_{i}=\underbrace{\prod_{i\in\left[
n\right]  \setminus\left\{  k\right\}  }}_{=\prod_{\substack{i\in\left[
n\right]  ;\\i\neq k}}}x_{i}=\prod_{\substack{i\in\left[  n\right]  ;\\i\neq
k}}x_{i}\\
&  =\prod_{\substack{j\in\left[  n\right]  ;\\j\neq k}}x_{j}%
\ \ \ \ \ \ \ \ \ \ \left(
\begin{array}
[c]{c}%
\text{here, we have renamed the index }i\text{ as }j\\
\text{in the product}%
\end{array}
\right) \\
&  =\prod_{\substack{j\in\left\{  1,2,\ldots,n\right\}  ;\\j\neq k}%
}x_{j}\ \ \ \ \ \ \ \ \ \ \left(  \text{since }\left[  n\right]  =\left\{
1,2,\ldots,n\right\}  \right)  .
\end{align*}
Comparing this with%
\[
y_{k}=\prod_{\substack{j\in\left\{  1,2,\ldots,n\right\}  ;\\j\neq k}%
}x_{j}\ \ \ \ \ \ \ \ \ \ \left(  \text{by the definition of }y_{k}\right)  ,
\]
we find $f_{n-1}=y_{k}$. Thus, (\ref{pf.prop.sol.vander-hook.viete-y.fn-1=})
is proven.]

For every $j\in\mathbb{N}$, we have%
\begin{equation}
e_{j}\left(  x_{1},x_{2},\ldots,x_{n}\right)  =\sum_{I\in\mathcal{P}%
_{j}\left(  \left[  n\right]  \right)  }\ \ \prod_{i\in I}x_{i}
\label{pf.prop.sol.vander-hook.viete-y.ej=}%
\end{equation}
\footnote{\textit{Proof of (\ref{pf.prop.sol.vander-hook.viete-y.ej=}):} Let
$j\in\mathbb{N}$. Recall that $\mathcal{P}_{j}\left(  \left[  n\right]
\right)  $ is the set of all $j$-element subsets of $\left[  n\right]  $ (by
the definition of $\mathcal{P}_{j}\left(  \left[  n\right]  \right)  $). In
other words, $\mathcal{P}_{j}\left(  \left[  n\right]  \right)  $ is the set
of all subsets $I$ of $\left[  n\right]  $ satisfying $\left\vert I\right\vert
=j$. Hence, $\sum_{I\in\mathcal{P}_{j}\left(  \left[  n\right]  \right)
}=\sum_{\substack{I\subseteq\left[  n\right]  ;\\\left\vert I\right\vert =j}}$
(an equality of summation signs). Thus,%
\[
\underbrace{\sum_{I\in\mathcal{P}_{j}\left(  \left[  n\right]  \right)  }%
}_{=\sum_{\substack{I\subseteq\left[  n\right]  ;\\\left\vert I\right\vert
=j}}}\ \ \prod_{i\in I}x_{i}=\sum_{\substack{I\subseteq\left[  n\right]
;\\\left\vert I\right\vert =j}}\ \ \prod_{i\in I}x_{i}.
\]
Comparing this with
\[
e_{j}\left(  x_{1},x_{2},\ldots,x_{n}\right)  =\sum_{\substack{I\subseteq
\left[  n\right]  ;\\\left\vert I\right\vert =j}}\ \ \prod_{i\in I}%
x_{i}\ \ \ \ \ \ \ \ \ \ \left(  \text{by the definition of }e_{j}\left(
x_{1},x_{2},\ldots,x_{n}\right)  \right)  ,
\]
we obtain $e_{j}\left(  x_{1},x_{2},\ldots,x_{n}\right)  =\sum_{I\in
\mathcal{P}_{j}\left(  \left[  n\right]  \right)  }\ \ \prod_{i\in I}x_{i}$.
This proves (\ref{pf.prop.sol.vander-hook.viete-y.ej=}).}. Also,%
\begin{equation}
e_{0}\left(  x_{1},x_{2},\ldots,x_{n}\right)  =1
\label{pf.prop.sol.vander-hook.viete-y.e0=}%
\end{equation}
\footnote{\textit{Proof of (\ref{pf.prop.sol.vander-hook.viete-y.e0=}):}
Proposition \ref{prop.sol.prop.binom.subsets.Pm.lem} \textbf{(a)} (applied to
$S=\left[  n\right]  $) yields $\mathcal{P}_{0}\left(  \left[  n\right]
\right)  =\left\{  \varnothing\right\}  $. Now,
(\ref{pf.prop.sol.vander-hook.viete-y.ej=}) (applied to $j=0$) yields%
\begin{align*}
e_{0}\left(  x_{1},x_{2},\ldots,x_{n}\right)   &  =\sum_{I\in\mathcal{P}%
_{0}\left(  \left[  n\right]  \right)  }\ \ \prod_{i\in I}x_{i}=\sum
_{I\in\left\{  \varnothing\right\}  }\ \ \prod_{i\in I}x_{i}%
\ \ \ \ \ \ \ \ \ \ \left(  \text{since }\mathcal{P}_{0}\left(  \left[
n\right]  \right)  =\left\{  \varnothing\right\}  \right) \\
&  =\prod_{i\in\varnothing}x_{i}=\left(  \text{empty product}\right)  =1.
\end{align*}
This proves (\ref{pf.prop.sol.vander-hook.viete-y.e0=}).}.

We shall now show that%
\begin{equation}
e_{m}\left(  x_{1},x_{2},\ldots,x_{n}\right)  =f_{m}+x_{k}f_{m-1}
\label{pf.prop.sol.vander-hook.viete-y.main}%
\end{equation}
for every positive integer $m$.

[\textit{Proof of (\ref{pf.prop.sol.vander-hook.viete-y.main}):} Let $m$ be a
positive integer. Thus, $m\in\mathbb{N}$ and $m-1\in\mathbb{N}$. The
definition of $f_{m}$ yields%
\begin{equation}
f_{m}=\sum_{I\in\mathcal{P}_{m}\left(  K\right)  }\ \ \prod_{i\in I}x_{i}.
\label{pf.prop.sol.vander-hook.viete-y.main.pf.fm=}%
\end{equation}
The definition of $f_{m-1}$ yields%
\begin{equation}
f_{m-1}=\sum_{I\in\mathcal{P}_{m-1}\left(  K\right)  }\ \ \prod_{i\in I}x_{i}.
\label{pf.prop.sol.vander-hook.viete-y.main.pf.fm-1=}%
\end{equation}

Recall that $k\in\left[  n\right]  $. Hence, we can apply Proposition
\ref{prop.sol.prop.binom.subsets.Pm.lem} \textbf{(c)} to $\left[  n\right]  $
and $k$ instead of $S$ and $s$. As a result, we obtain the following two results:

\begin{itemize}
\item We have $\mathcal{P}_{m}\left(  \left[  n\right]  \setminus\left\{
k\right\}  \right)  \subseteq\mathcal{P}_{m}\left(  \left[  n\right]  \right)
$.

\item The map%
\begin{align*}
\mathcal{P}_{m-1}\left(  \left[  n\right]  \setminus\left\{  k\right\}
\right)   &  \rightarrow\mathcal{P}_{m}\left(  \left[  n\right]  \right)
\setminus\mathcal{P}_{m}\left(  \left[  n\right]  \setminus\left\{  k\right\}
\right)  ,\\
U  &  \mapsto U\cup\left\{  k\right\}
\end{align*}
is well-defined and a bijection.
\end{itemize}

Since $\left[  n\right]  \setminus\left\{  k\right\}  =K$, these two results
can be rewritten as follows:

\begin{itemize}
\item We have $\mathcal{P}_{m}\left(  K\right)  \subseteq\mathcal{P}%
_{m}\left(  \left[  n\right]  \right)  $.

\item The map%
\begin{align*}
\mathcal{P}_{m-1}\left(  K\right)   &  \rightarrow\mathcal{P}_{m}\left(
\left[  n\right]  \right)  \setminus\mathcal{P}_{m}\left(  K\right)  ,\\
U  &  \mapsto U\cup\left\{  k\right\}
\end{align*}
is well-defined and a bijection.
\end{itemize}

Furthermore, every $I\in\mathcal{P}_{m-1}\left(  K\right)  $ satisfies%
\begin{equation}
\prod_{i\in I\cup\left\{  k\right\}  }x_{i}=x_{k}\prod_{i\in I}x_{i}
\label{pf.prop.sol.vander-hook.viete-y.main.pf.2}%
\end{equation}
\footnote{\textit{Proof of (\ref{pf.prop.sol.vander-hook.viete-y.main.pf.2}):}
Let $I\in\mathcal{P}_{m-1}\left(  K\right)  $. Thus, $I$ is an $\left(
m-1\right)  $-element subset of $K$ (since $\mathcal{P}_{m-1}\left(  K\right)
$ is the set of all $\left(  m-1\right)  $-element subsets of $K$ (by the
definition of $\mathcal{P}_{m-1}\left(  K\right)  $)). In particular, $I$ is a
subset of $K$. Thus, $I\subseteq K=\left[  n\right]  \setminus\left\{
k\right\}  $.
\par
We have $k\in\left\{  k\right\}  $ and thus $k\notin\left[  n\right]
\setminus\left\{  k\right\}  $. If we had $k\in I$, then we would have $k\in
I\subseteq\left[  n\right]  \setminus\left\{  k\right\}  $, which would
contradict $k\notin\left[  n\right]  \setminus\left\{  k\right\}  $. Hence, we
cannot have $k\in I$. Thus, we have $k\notin I$. Hence, the sets $I$ and
$\left\{  k\right\}  $ are disjoint. Thus,%
\[
\prod_{i\in I\cup\left\{  k\right\}  }x_{i}=\left(  \prod_{i\in I}%
x_{i}\right)  \underbrace{\left(  \prod_{i\in\left\{  k\right\}  }%
x_{i}\right)  }_{=x_{k}}=\left(  \prod_{i\in I}x_{i}\right)  x_{k}=x_{k}%
\prod_{i\in I}x_{i}.
\]
This proves (\ref{pf.prop.sol.vander-hook.viete-y.main.pf.2}).}. Now,%
\begin{align}
&  \sum_{I\in\mathcal{P}_{m}\left(  \left[  n\right]  \right)  \setminus
\mathcal{P}_{m}\left(  K\right)  }\ \ \prod_{i\in I}x_{i}\nonumber\\
&  =\sum_{U\in\mathcal{P}_{m-1}\left(  K\right)  }\ \ \prod_{i\in
U\cup\left\{  k\right\}  }x_{i}\nonumber\\
&  \ \ \ \ \ \ \ \ \ \ \ \ \ \ \ \ \ \ \ \ \left(
\begin{array}
[c]{c}%
\text{here, we have substituted }U\cup\left\{  k\right\}  \text{ for }I\text{
in the sum, since}\\
\text{the map }\mathcal{P}_{m-1}\left(  K\right)  \rightarrow\mathcal{P}%
_{m}\left(  \left[  n\right]  \right)  \setminus\mathcal{P}_{m}\left(
K\right)  ,\ U\mapsto U\cup\left\{  k\right\} \\
\text{is a bijection}%
\end{array}
\right) \nonumber\\
&  =\sum_{I\in\mathcal{P}_{m-1}\left(  K\right)  }\ \ \underbrace{\prod_{i\in
I\cup\left\{  k\right\}  }x_{i}}_{\substack{=x_{k}\prod_{i\in I}%
x_{i}\\\text{(by (\ref{pf.prop.sol.vander-hook.viete-y.main.pf.2}))}%
}}\ \ \ \ \ \ \ \ \ \ \left(
\begin{array}
[c]{c}%
\text{here, we have renamed the}\\
\text{summation index }U\text{ as }I
\end{array}
\right) \nonumber\\
&  =\sum_{I\in\mathcal{P}_{m-1}\left(  K\right)  }x_{k}\prod_{i\in I}%
x_{i}=x_{k}\sum_{I\in\mathcal{P}_{m-1}\left(  K\right)  }\ \ \prod_{i\in
I}x_{i}. \label{pf.prop.sol.vander-hook.viete-y.main.pf.2b}%
\end{align}

From (\ref{pf.prop.sol.vander-hook.viete-y.ej=}) (applied to $j=m$), we obtain%
\begin{align*}
e_{m}\left(  x_{1},x_{2},\ldots,x_{n}\right)   &  =\sum_{I\in\mathcal{P}%
_{m}\left(  \left[  n\right]  \right)  }\ \ \prod_{i\in I}x_{i}%
=\underbrace{\sum_{\substack{I\in\mathcal{P}_{m}\left(  \left[  n\right]
\right)  ;\\I\in\mathcal{P}_{m}\left(  K\right)  }}}_{\substack{=\sum
_{I\in\mathcal{P}_{m}\left(  K\right)  }\\\text{(since }\mathcal{P}_{m}\left(
K\right)  \subseteq\mathcal{P}_{m}\left(  \left[  n\right]  \right)  \text{)}%
}}\prod_{i\in I}x_{i}+\underbrace{\sum_{\substack{I\in\mathcal{P}_{m}\left(
\left[  n\right]  \right)  ;\\I\notin\mathcal{P}_{m}\left(  K\right)  }%
}}_{=\sum_{I\in\mathcal{P}_{m}\left(  \left[  n\right]  \right)
\setminus\mathcal{P}_{m}\left(  K\right)  }}\prod_{i\in I}x_{i}\\
&  \ \ \ \ \ \ \ \ \ \ \ \ \ \ \ \ \ \ \ \ \left(
\begin{array}
[c]{c}%
\text{since each }I\in\mathcal{P}_{m}\left(  \left[  n\right]  \right)  \text{
satisfies either }I\in\mathcal{P}_{m}\left(  K\right) \\
\text{or }I\notin\mathcal{P}_{m}\left(  K\right)  \text{ (but not both)}%
\end{array}
\right) \\
&  =\sum_{I\in\mathcal{P}_{m}\left(  K\right)  }\ \ \prod_{i\in I}%
x_{i}+\underbrace{\sum_{I\in\mathcal{P}_{m}\left(  \left[  n\right]  \right)
\setminus\mathcal{P}_{m}\left(  K\right)  }\ \ \prod_{i\in I}x_{i}%
}_{\substack{=x_{k}\sum_{I\in\mathcal{P}_{m-1}\left(  K\right)  }%
\ \ \prod_{i\in I}x_{i}\\\text{(by
(\ref{pf.prop.sol.vander-hook.viete-y.main.pf.2b}))}}}\\
&  =\underbrace{\sum_{I\in\mathcal{P}_{m}\left(  K\right)  }\ \ \prod_{i\in
I}x_{i}}_{\substack{=f_{m}\\\text{(by
(\ref{pf.prop.sol.vander-hook.viete-y.main.pf.fm=}))}}}+x_{k}\underbrace{\sum
_{I\in\mathcal{P}_{m-1}\left(  K\right)  }\ \ \prod_{i\in I}x_{i}%
}_{\substack{=f_{m-1}\\\text{(by
(\ref{pf.prop.sol.vander-hook.viete-y.main.pf.fm-1=}))}}}\\
&  =f_{m}+x_{k}f_{m-1}.
\end{align*}
This proves (\ref{pf.prop.sol.vander-hook.viete-y.main}).]

Recall that $k\in\left\{  1,2,\ldots,n\right\}  $. Hence, $1\leq k\leq n$, so
that $1\leq n$. Thus, $n\geq1$, so that $n-1\in\mathbb{N}$. Hence,
$0\in\left\{  0,1,\ldots,n-1\right\}  $ and $n-1\in\left\{  0,1,\ldots
,n-1\right\}  $.

Now,%
\begin{align}
&  \sum_{j=0}^{n-1}\left(  -1\right)  ^{n-1-j}e_{j}\left(  x_{1},x_{2}%
,\ldots,x_{n}\right)  x_{k}^{n-1-j}\nonumber\\
&  =\underbrace{\left(  -1\right)  ^{n-1-0}}_{=\left(  -1\right)  ^{n-1}%
}\underbrace{e_{0}\left(  x_{1},x_{2},\ldots,x_{n}\right)  }%
_{\substack{=1\\\text{(by (\ref{pf.prop.sol.vander-hook.viete-y.e0=}))}%
}}\underbrace{x_{k}^{n-1-0}}_{=x_{k}^{n-1}}+\sum_{j=1}^{n-1}\left(  -1\right)
^{n-1-j}\underbrace{e_{j}\left(  x_{1},x_{2},\ldots,x_{n}\right)
}_{\substack{=f_{j}+x_{k}f_{j-1}\\\text{(by
(\ref{pf.prop.sol.vander-hook.viete-y.main}) (applied to }m=j\text{))}}%
}x_{k}^{n-1-j}\nonumber\\
&  \ \ \ \ \ \ \ \ \ \ \ \ \ \ \ \ \ \ \ \ \left(
\begin{array}
[c]{c}%
\text{here, we have split off the addend for }j=0\text{ from the sum,}\\
\text{since }0\in\left\{  0,1,\ldots,n-1\right\}
\end{array}
\right) \nonumber\\
&  =\left(  -1\right)  ^{n-1}x_{k}^{n-1}+\sum_{j=1}^{n-1}\underbrace{\left(
-1\right)  ^{n-1-j}\left(  f_{j}+x_{k}f_{j-1}\right)  x_{k}^{n-1-j}}_{=\left(
-1\right)  ^{n-1-j}f_{j}x_{k}^{n-1-j}+\left(  -1\right)  ^{n-1-j}x_{k}%
f_{j-1}x_{k}^{n-1-j}}\nonumber\\
&  =\left(  -1\right)  ^{n-1}x_{k}^{n-1}+\underbrace{\sum_{j=1}^{n-1}\left(
\left(  -1\right)  ^{n-1-j}f_{j}x_{k}^{n-1-j}+\left(  -1\right)  ^{n-1-j}%
x_{k}f_{j-1}x_{k}^{n-1-j}\right)  }_{=\sum_{j=1}^{n-1}\left(  -1\right)
^{n-1-j}f_{j}x_{k}^{n-1-j}+\sum_{j=1}^{n-1}\left(  -1\right)  ^{n-1-j}%
x_{k}f_{j-1}x_{k}^{n-1-j}}\nonumber\\
&  =\left(  -1\right)  ^{n-1}x_{k}^{n-1}+\sum_{j=1}^{n-1}\left(  -1\right)
^{n-1-j}f_{j}x_{k}^{n-1-j}+\sum_{j=1}^{n-1}\left(  -1\right)  ^{n-1-j}%
x_{k}f_{j-1}x_{k}^{n-1-j}. \label{pf.prop.sol.vander-hook.viete-y.5}%
\end{align}

But
\begin{align}
&  \sum_{j=1}^{n-1}\left(  -1\right)  ^{n-1-j}\underbrace{x_{k}f_{j-1}%
}_{=f_{j-1}x_{k}}x_{k}^{n-1-j}\nonumber\\
&  =\sum_{j=1}^{n-1}\left(  -1\right)  ^{n-1-j}f_{j-1}\underbrace{x_{k}%
x_{k}^{n-1-j}}_{\substack{=x_{k}^{\left(  n-1-j\right)  +1}=x_{k}%
^{n-j}\\\text{(since }\left(  n-1-j\right)  +1=n-j\text{)}}}=\sum_{j=1}%
^{n-1}\left(  -1\right)  ^{n-1-j}f_{j-1}x_{k}^{n-j}\nonumber\\
&  =\sum_{j=0}^{\left(  n-1\right)  -1}\underbrace{\left(  -1\right)
^{n-1-\left(  j+1\right)  }}_{\substack{=\left(  -1\right)  ^{n-j}%
\\\text{(since }n-1-\left(  j+1\right)  =n-j-2\equiv n-j\operatorname{mod}%
2\text{)}}}\underbrace{f_{\left(  j+1\right)  -1}}_{\substack{=f_{j}%
\\\text{(since }\left(  j+1\right)  -1=j\text{)}}}\underbrace{x_{k}^{n-\left(
j+1\right)  }}_{\substack{=x_{k}^{n-1-j}\\\text{(since }n-\left(  j+1\right)
=n-1-j\text{)}}}\nonumber\\
&  \ \ \ \ \ \ \ \ \ \ \ \ \ \ \ \ \ \ \ \ \left(  \text{here, we have
substituted }j+1\text{ for }j\text{ in the sum}\right) \nonumber\\
&  =\sum_{j=0}^{\left(  n-1\right)  -1}\underbrace{\left(  -1\right)  ^{n-j}%
}_{\substack{=\left(  -1\right)  \left(  -1\right)  ^{n-j-1}\\=-\left(
-1\right)  ^{n-j-1}=-\left(  -1\right)  ^{n-1-j}\\\text{(since }%
n-j-1=n-1-j\text{)}}}f_{j}x_{k}^{n-1-j}=\sum_{j=0}^{\left(  n-1\right)
-1}\left(  -\left(  -1\right)  ^{n-1-j}\right)  f_{j}x_{k}^{n-1-j}\nonumber\\
&  =-\sum_{j=0}^{\left(  n-1\right)  -1}\left(  -1\right)  ^{n-1-j}f_{j}%
x_{k}^{n-1-j}. \label{pf.prop.sol.vander-hook.viete-y.6a}%
\end{align}
Also,%
\begin{align*}
&  \sum_{j=0}^{n-1}\left(  -1\right)  ^{n-1-j}f_{j}x_{k}^{n-1-j}\\
&  =\underbrace{\left(  -1\right)  ^{n-1-0}}_{=\left(  -1\right)  ^{n-1}%
}\underbrace{f_{0}}_{\substack{=1\\\text{(by
(\ref{pf.prop.sol.vander-hook.viete-y.f0=}))}}}\underbrace{x_{k}^{n-1-0}%
}_{=x_{k}^{n-1}}+\sum_{j=1}^{n-1}\left(  -1\right)  ^{n-1-j}f_{j}x_{k}%
^{n-1-j}\\
&  \ \ \ \ \ \ \ \ \ \ \ \ \ \ \ \ \ \ \ \ \left(
\begin{array}
[c]{c}%
\text{here, we have split off the addend for }j=0\text{ from the sum,}\\
\text{since }0\in\left\{  0,1,\ldots,n-1\right\}
\end{array}
\right) \\
&  =\left(  -1\right)  ^{n-1}x_{k}^{n-1}+\sum_{j=1}^{n-1}\left(  -1\right)
^{n-1-j}f_{j}x_{k}^{n-1-j},
\end{align*}
so that%
\begin{align}
&  \left(  -1\right)  ^{n-1}x_{k}^{n-1}+\sum_{j=1}^{n-1}\left(  -1\right)
^{n-1-j}f_{j}x_{k}^{n-1-j}\nonumber\\
&  =\sum_{j=0}^{n-1}\left(  -1\right)  ^{n-1-j}f_{j}x_{k}^{n-1-j}\nonumber\\
&  =\sum_{j=0}^{\left(  n-1\right)  -1}\left(  -1\right)  ^{n-1-j}f_{j}%
x_{k}^{n-1-j}+\underbrace{\left(  -1\right)  ^{n-1-\left(  n-1\right)  }%
}_{\substack{=\left(  -1\right)  ^{0}\\\text{(since }n-1-\left(  n-1\right)
=0\text{)}}}\underbrace{f_{n-1}}_{\substack{=y_{k}\\\text{(by
(\ref{pf.prop.sol.vander-hook.viete-y.fn-1=}))}}}\underbrace{x_{k}%
^{n-1-\left(  n-1\right)  }}_{\substack{=x_{k}^{0}\\\text{(since }n-1-\left(
n-1\right)  =0\text{)}}}\nonumber\\
&  \ \ \ \ \ \ \ \ \ \ \ \ \ \ \ \ \ \ \ \ \left(
\begin{array}
[c]{c}%
\text{here, we have split off the addend for }j=n-1\text{ from the sum,}\\
\text{since }n-1\in\left\{  0,1,\ldots,n-1\right\}
\end{array}
\right) \nonumber\\
&  =\sum_{j=0}^{\left(  n-1\right)  -1}\left(  -1\right)  ^{n-1-j}f_{j}%
x_{k}^{n-1-j}+\underbrace{\left(  -1\right)  ^{0}}_{=1}y_{k}\underbrace{x_{k}%
^{0}}_{=1}\nonumber\\
&  =\sum_{j=0}^{\left(  n-1\right)  -1}\left(  -1\right)  ^{n-1-j}f_{j}%
x_{k}^{n-1-j}+y_{k}. \label{pf.prop.sol.vander-hook.viete-y.6b}%
\end{align}

Now, (\ref{pf.prop.sol.vander-hook.viete-y.5}) becomes%
\begin{align*}
&  \sum_{j=0}^{n-1}\left(  -1\right)  ^{n-1-j}e_{j}\left(  x_{1},x_{2}%
,\ldots,x_{n}\right)  x_{k}^{n-1-j}\\
&  =\underbrace{\left(  -1\right)  ^{n-1}x_{k}^{n-1}+\sum_{j=1}^{n-1}\left(
-1\right)  ^{n-1-j}f_{j}x_{k}^{n-1-j}}_{\substack{=\sum_{j=0}^{\left(
n-1\right)  -1}\left(  -1\right)  ^{n-1-j}f_{j}x_{k}^{n-1-j}+y_{k}\\\text{(by
(\ref{pf.prop.sol.vander-hook.viete-y.6b}))}}}+\underbrace{\sum_{j=1}%
^{n-1}\left(  -1\right)  ^{n-1-j}x_{k}f_{j-1}x_{k}^{n-1-j}}_{\substack{=-\sum
_{j=0}^{\left(  n-1\right)  -1}\left(  -1\right)  ^{n-1-j}f_{j}x_{k}%
^{n-1-j}\\\text{(by (\ref{pf.prop.sol.vander-hook.viete-y.6a}))}}}\\
&  =\sum_{j=0}^{\left(  n-1\right)  -1}\left(  -1\right)  ^{n-1-j}f_{j}%
x_{k}^{n-1-j}+y_{k}+\left(  -\sum_{j=0}^{\left(  n-1\right)  -1}\left(
-1\right)  ^{n-1-j}f_{j}x_{k}^{n-1-j}\right) \\
&  =y_{k}.
\end{align*}
This proves Proposition \ref{prop.sol.vander-hook.viete-y}.
\end{proof}
\end{verlong}

\begin{lemma}
\label{lem.sol.vander-hook.y-sum.1}Let $n\in\mathbb{N}$. Let $x_{1}%
,x_{2},\ldots,x_{n}$ be $n$ elements of $\mathbb{K}$. For each $i\in\left\{
1,2,\ldots,n\right\}  $, set $y_{i}=\prod_{\substack{j\in\left\{
1,2,\ldots,n\right\}  ;\\j\neq i}}x_{j}$.

Let $A$ be the $n\times n$-matrix $\left(  x_{i}^{n-j}\right)  _{1\leq i\leq
n,\ 1\leq j\leq n}$. Let $q\in\left\{  1,2,\ldots,n\right\}  $. For every
$j\in\mathbb{N}$, define an element $\mathfrak{e}_{j}\in\mathbb{K}$ by
$\mathfrak{e}_{j}=e_{j}\left(  x_{1},x_{2},\ldots,x_{n}\right)  $. Then,%
\[
\sum_{k=1}^{n}\left(  -1\right)  ^{k+q}y_{k}\det\left(  A_{\sim k,\sim
q}\right)  =\left(  -1\right)  ^{n-q}\mathfrak{e}_{q-1}\det A.
\]

\end{lemma}

\begin{proof}
[Proof of Lemma \ref{lem.sol.vander-hook.y-sum.1}.]We have $A=\left(
x_{i}^{n-j}\right)  _{1\leq i\leq n,\ 1\leq j\leq n}$ (by the definition of
$A$). Hence, Lemma \ref{lem.sol.vander-hook.lap0} (applied to $a_{i,j}%
=x_{i}^{n-j}$) yields that%
\begin{equation}
\sum_{p=1}^{n}\left(  -1\right)  ^{p+q}x_{p}^{n-r}\det\left(  A_{\sim p,\sim
q}\right)  =\delta_{q,r}\det A \label{pf.lem.sol.vander-hook.y-sum.1.0}%
\end{equation}
for each $r\in\left\{  1,2,\ldots,n\right\}  $.

Let $j\in\left\{  0,1,\ldots,n-1\right\}  $. Thus, $j+1\in\left\{
1,2,\ldots,n\right\}  $. Hence, (\ref{pf.lem.sol.vander-hook.y-sum.1.0})
(applied to $r=j+1$) yields%
\[
\sum_{p=1}^{n}\left(  -1\right)  ^{p+q}x_{p}^{n-\left(  j+1\right)  }%
\det\left(  A_{\sim p,\sim q}\right)  =\delta_{q,j+1}\det A.
\]
Hence,%
\begin{align}
\delta_{q,j+1}\det A  &  =\sum_{p=1}^{n}\left(  -1\right)  ^{p+q}%
\underbrace{x_{p}^{n-\left(  j+1\right)  }}_{\substack{=x_{p}^{n-1-j}%
\\\text{(since }n-\left(  j+1\right)  =n-1-j\text{)}}}\det\left(  A_{\sim
p,\sim q}\right) \nonumber\\
&  =\sum_{p=1}^{n}\left(  -1\right)  ^{p+q}x_{p}^{n-1-j}\det\left(  A_{\sim
p,\sim q}\right) \nonumber\\
&  =\sum_{k=1}^{n}\left(  -1\right)  ^{k+q}x_{k}^{n-1-j}\det\left(  A_{\sim
k,\sim q}\right)  \label{pf.lem.sol.vander-hook.y-sum.1.1}%
\end{align}
(here, we have renamed the summation index $p$ as $k$).

Now, forget that we fixed $j$. We thus have proven
(\ref{pf.lem.sol.vander-hook.y-sum.1.1}) for every $j\in\left\{
0,1,\ldots,n-1\right\}  $.

Every $k\in\mathbb{N}$ satisfies%
\begin{align*}
y_{k}  &  =\sum_{j=0}^{n-1}\left(  -1\right)  ^{n-1-j}\underbrace{e_{j}\left(
x_{1},x_{2},\ldots,x_{n}\right)  }_{\substack{=\mathfrak{e}_{j}\\\text{(since
}\mathfrak{e}_{j}=e_{j}\left(  x_{1},x_{2},\ldots,x_{n}\right)  \\\text{(by
the definition of }\mathfrak{e}_{j}\text{))}}}x_{k}^{n-1-j}%
\ \ \ \ \ \ \ \ \ \ \left(  \text{by Proposition
\ref{prop.sol.vander-hook.viete-y}}\right) \\
&  =\sum_{j=0}^{n-1}\left(  -1\right)  ^{n-1-j}\mathfrak{e}_{j}x_{k}^{n-1-j}.
\end{align*}
Thus,%
\begin{align*}
&  \sum_{k=1}^{n}\left(  -1\right)  ^{k+q}\underbrace{y_{k}}_{=\sum
_{j=0}^{n-1}\left(  -1\right)  ^{n-1-j}\mathfrak{e}_{j}x_{k}^{n-1-j}}%
\det\left(  A_{\sim k,\sim q}\right) \\
&  =\sum_{k=1}^{n}\left(  -1\right)  ^{k+q}\left(  \sum_{j=0}^{n-1}\left(
-1\right)  ^{n-1-j}\mathfrak{e}_{j}x_{k}^{n-1-j}\right)  \det\left(  A_{\sim
k,\sim q}\right) \\
&  =\underbrace{\sum_{k=1}^{n}\ \ \sum_{j=0}^{n-1}}_{=\sum_{j=0}^{n-1}%
\ \ \sum_{k=1}^{n}}\left(  -1\right)  ^{k+q}\left(  -1\right)  ^{n-1-j}%
\mathfrak{e}_{j}x_{k}^{n-1-j}\det\left(  A_{\sim k,\sim q}\right) \\
&  =\sum_{j=0}^{n-1}\ \ \underbrace{\sum_{k=1}^{n}\left(  -1\right)
^{k+q}\left(  -1\right)  ^{n-1-j}\mathfrak{e}_{j}x_{k}^{n-1-j}\det\left(
A_{\sim k,\sim q}\right)  }_{=\left(  -1\right)  ^{n-1-j}\mathfrak{e}_{j}%
\sum_{k=1}^{n}\left(  -1\right)  ^{k+q}x_{k}^{n-1-j}\det\left(  A_{\sim k,\sim
q}\right)  }\\
&  =\sum_{j=0}^{n-1}\left(  -1\right)  ^{n-1-j}\mathfrak{e}_{j}%
\underbrace{\sum_{k=1}^{n}\left(  -1\right)  ^{k+q}x_{k}^{n-1-j}\det\left(
A_{\sim k,\sim q}\right)  }_{\substack{=\delta_{q,j+1}\det A\\\text{(by
(\ref{pf.lem.sol.vander-hook.y-sum.1.1}))}}}\\
&  =\sum_{j=0}^{n-1}\left(  -1\right)  ^{n-1-j}\mathfrak{e}_{j}\delta
_{q,j+1}\det A=\underbrace{\sum_{j=1}^{n}}_{=\sum_{j\in\left\{  1,2,\ldots
,n\right\}  }}\underbrace{\left(  -1\right)  ^{n-1-\left(  j-1\right)  }%
}_{\substack{=\left(  -1\right)  ^{n-j}\\\text{(since }n-1-\left(  j-1\right)
=n-j\text{)}}}\mathfrak{e}_{j-1}\underbrace{\delta_{q,\left(  j-1\right)  +1}%
}_{\substack{=\delta_{q,j}\\\text{(since }\left(  j-1\right)  +1=j\text{)}%
}}\det A\\
&  \ \ \ \ \ \ \ \ \ \ \left(  \text{here, we have substituted }j-1\text{ for
}j\text{ in the sum}\right) \\
&  =\sum_{j\in\left\{  1,2,\ldots,n\right\}  }\left(  -1\right)
^{n-j}\mathfrak{e}_{j-1}\delta_{q,j}\det A
\end{align*}%
\begin{align*}
&  =\left(  -1\right)  ^{n-q}\mathfrak{e}_{q-1}\underbrace{\delta_{q,q}%
}_{\substack{=1\\\text{(since }q=q\text{)}}}\det A+\sum_{\substack{j\in
\left\{  1,2,\ldots,n\right\}  ;\\j\neq q}}\left(  -1\right)  ^{n-j}%
\mathfrak{e}_{j-1}\underbrace{\delta_{q,j}}_{\substack{=0\\\text{(since }q\neq
j\\\text{(since }j\neq q\text{))}}}\det A\\
&  \ \ \ \ \ \ \ \ \ \ \left(
\begin{array}
[c]{c}%
\text{here, we have split off the addend for }j=q\text{ from the sum}\\
\text{(since }q\in\left\{  1,2,\ldots,n\right\}  \text{)}%
\end{array}
\right) \\
&  =\left(  -1\right)  ^{n-q}\mathfrak{e}_{q-1}\det A+\underbrace{\sum
_{\substack{j\in\left\{  1,2,\ldots,n\right\}  ;\\j\neq q}}\left(  -1\right)
^{n-j}\mathfrak{e}_{j-1}0\det A}_{=0}=\left(  -1\right)  ^{n-q}\mathfrak{e}%
_{q-1}\det A.
\end{align*}
This proves Lemma \ref{lem.sol.vander-hook.y-sum.1}.
\end{proof}

\begin{lemma}
\label{lem.sol.vander-hook.y-sum.2}Let $n\in\mathbb{N}$. Let $x_{1}%
,x_{2},\ldots,x_{n}$ be $n$ elements of $\mathbb{K}$. For each $i\in\left\{
1,2,\ldots,n\right\}  $, set $y_{i}=\prod_{\substack{j\in\left\{
1,2,\ldots,n\right\}  ;\\j\neq i}}x_{j}$.

Let $A$ be the $n\times n$-matrix $\left(  x_{i}^{n-j}\right)  _{1\leq i\leq
n,\ 1\leq j\leq n}$. Let $q\in\left\{  1,2,\ldots,n\right\}  $ and $\ell
\in\left\{  0,1,\ldots,n-q\right\}  $. For every $j\in\mathbb{N}$, define an
element $\mathfrak{e}_{j}\in\mathbb{K}$ by $\mathfrak{e}_{j}=e_{j}\left(
x_{1},x_{2},\ldots,x_{n}\right)  $. Then,%
\[
\sum_{k=1}^{n}\left(  -1\right)  ^{k+q}y_{k}x_{k}^{n-q-\ell}\det\left(
A_{\sim k,\sim q}\right)  =\delta_{\ell,n-q}\left(  -1\right)  ^{n-q}%
\mathfrak{e}_{q-1}\det A.
\]

\end{lemma}

\begin{proof}
[Proof of Lemma \ref{lem.sol.vander-hook.y-sum.2}.]We are in one of the
following two cases:

\textit{Case 1:} We have $\ell=n-q$.

\textit{Case 2:} We have $\ell\neq n-q$.

Let us first consider Case 1. In this case, we have $\ell=n-q$. Thus,
$n-q-\underbrace{\ell}_{=n-q}=n-q-\left(  n-q\right)  =0$. Now,%
\begin{align*}
&  \sum_{k=1}^{n}\left(  -1\right)  ^{k+q}y_{k}\underbrace{x_{k}^{n-q-\ell}%
}_{\substack{=x_{k}^{0}\\\text{(since }n-q-\ell=0\text{)}}}\det\left(  A_{\sim
k,\sim q}\right) \\
&  =\sum_{k=1}^{n}\left(  -1\right)  ^{k+q}y_{k}\underbrace{x_{k}^{0}}%
_{=1}\det\left(  A_{\sim k,\sim q}\right) \\
&  =\sum_{k=1}^{n}\left(  -1\right)  ^{k+q}y_{k}\det\left(  A_{\sim k,\sim
q}\right)  =\left(  -1\right)  ^{n-q}\mathfrak{e}_{q-1}\det A
\end{align*}
(by Lemma \ref{lem.sol.vander-hook.y-sum.1}). Comparing this with%
\[
\underbrace{\delta_{\ell,n-q}}_{\substack{=1\\\text{(since }\ell=n-q\text{)}%
}}\left(  -1\right)  ^{n-q}\mathfrak{e}_{q-1}\det A=\left(  -1\right)
^{n-q}\mathfrak{e}_{q-1}\det A,
\]
we obtain%
\[
\sum_{k=1}^{n}\left(  -1\right)  ^{k+q}y_{k}x_{k}^{n-q-\ell}\det\left(
A_{\sim k,\sim q}\right)  =\delta_{\ell,n-q}\left(  -1\right)  ^{n-q}%
\mathfrak{e}_{q-1}\det A.
\]
Hence, Lemma \ref{lem.sol.vander-hook.y-sum.2} is proven in Case 1.

Let us now consider Case 2. In this case, we have $\ell\neq n-q$. Combining
this with $\ell\leq n-q$ (since $\ell\in\left\{  0,1,\ldots,n-q\right\}  $),
we obtain $\ell<n-q$. Hence, $\ell\leq\left(  n-q\right)  -1$ (since $\ell$
and $n-q$ are integers). In other words, $\ell+1\leq n-q$. But $\ell\geq0$
(since $\ell\in\left\{  0,1,\ldots,n-q\right\}  $), so that $\ell+1\geq1$.
Combining this with $\ell+1\leq n-q$, we obtain $\ell+1\in\left\{
1,2,\ldots,n-q\right\}  $. Thus, Lemma \ref{lem.sol.vander-hook.lap2.variant}
(applied to $\ell+1$ instead of $\ell$) yields%
\begin{equation}
\sum_{k=1}^{n}\left(  -1\right)  ^{k+q}x_{k}^{n-q-\left(  \ell+1\right)  }%
\det\left(  A_{\sim k,\sim q}\right)  =0.
\label{pf.lem.sol.vander-hook.y-sum.2.c2.1}%
\end{equation}

But every $k\in\left\{  1,2,\ldots,n\right\}  $ satisfies%
\begin{equation}
y_{k}x_{k}^{n-q-\ell}=\mathfrak{e}_{n}x_{k}^{n-q-\left(  \ell+1\right)  }
\label{pf.lem.sol.vander-hook.y-sum.2.c2.3}%
\end{equation}
\footnote{\textit{Proof of (\ref{pf.lem.sol.vander-hook.y-sum.2.c2.3}):} Let
$k\in\left\{  1,2,\ldots,n\right\}  $. We have $n-q-\underbrace{\ell}%
_{\leq\left(  n-q\right)  -1}\geq n-q-\left(  \left(  n-q\right)  -1\right)
=1$. Thus, $x_{k}^{n-q-\ell}=x_{k}x_{k}^{\left(  n-q-\ell\right)  -1}%
=x_{k}x_{k}^{n-q-\left(  \ell+1\right)  }$ (since $\left(  n-q-\ell\right)
-1=n-q-\left(  \ell+1\right)  $).
\par
The definition of $y_{k}$ yields $y_{k}=\prod_{\substack{j\in\left\{
1,2,\ldots,n\right\}  ;\\j\neq k}}x_{j}$. But the definition of $\mathfrak{e}%
_{n}$ yields $\mathfrak{e}_{n}=e_{n}\left(  x_{1},x_{2},\ldots,x_{n}\right)
=\prod_{i=1}^{n}x_{i}$ (by Proposition \ref{prop.sol.vander-hook.viete}
\textbf{(b)} (applied to $t=0$)). Hence,%
\begin{align}
\mathfrak{e}_{n}  &  =\prod_{i=1}^{n}x_{i}=\underbrace{\prod_{j=1}^{n}%
}_{=\prod_{j\in\left\{  1,2,\ldots,n\right\}  }}x_{j}%
\ \ \ \ \ \ \ \ \ \ \left(
\begin{array}
[c]{c}%
\text{here, we have renamed the index }i\text{ as }j\\
\text{in the product}%
\end{array}
\right) \nonumber\\
&  =\prod_{j\in\left\{  1,2,\ldots,n\right\}  }x_{j}=x_{k}\underbrace{\prod
_{\substack{j\in\left\{  1,2,\ldots,n\right\}  ;\\j\neq k}}x_{j}}_{=y_{k}%
}\ \ \ \ \ \ \ \ \ \ \left(
\begin{array}
[c]{c}%
\text{here, we have split off the factor}\\
\text{for }j=k\text{ from the product}%
\end{array}
\right) \nonumber\\
&  =x_{k}y_{k}=y_{k}x_{k}. \label{pf.lem.sol.vander-hook.y-sum.2.c2.3.pf.2}%
\end{align}
Now,%
\[
y_{k}\underbrace{x_{k}^{n-q-\ell}}_{=x_{k}x_{k}^{n-q-\left(  \ell+1\right)  }%
}=\underbrace{y_{k}x_{k}}_{\substack{=\mathfrak{e}_{n}\\\text{(by
(\ref{pf.lem.sol.vander-hook.y-sum.2.c2.3.pf.2}))}}}x_{k}^{n-q-\left(
\ell+1\right)  }=\mathfrak{e}_{n}x_{k}^{n-q-\left(  \ell+1\right)  }.
\]
This proves (\ref{pf.lem.sol.vander-hook.y-sum.2.c2.3}).}. Now,%
\begin{align*}
&  \sum_{k=1}^{n}\left(  -1\right)  ^{k+q}\underbrace{y_{k}x_{k}^{n-q-\ell}%
}_{\substack{=\mathfrak{e}_{n}x_{k}^{n-q-\left(  \ell+1\right)  }\\\text{(by
(\ref{pf.lem.sol.vander-hook.y-sum.2.c2.3}))}}}\det\left(  A_{\sim k,\sim
q}\right) \\
&  =\sum_{k=1}^{n}\left(  -1\right)  ^{k+q}\mathfrak{e}_{n}x_{k}^{n-q-\left(
\ell+1\right)  }\det\left(  A_{\sim k,\sim q}\right) \\
&  =\mathfrak{e}_{n}\underbrace{\sum_{k=1}^{n}\left(  -1\right)  ^{k+q}%
x_{k}^{n-q-\left(  \ell+1\right)  }\det\left(  A_{\sim k,\sim q}\right)
}_{\substack{=0\\\text{(by (\ref{pf.lem.sol.vander-hook.y-sum.2.c2.1}))}%
}}=\mathfrak{e}_{n}0=0.
\end{align*}
Comparing this with%
\[
\underbrace{\delta_{\ell,n-q}}_{\substack{=0\\\text{(since }\ell\neq
n-q\text{)}}}\left(  -1\right)  ^{n-q}\mathfrak{e}_{q-1}\det A=0\left(
-1\right)  ^{n-q}\mathfrak{e}_{q-1}\det A=0,
\]
we obtain%
\[
\sum_{k=1}^{n}\left(  -1\right)  ^{k+q}y_{k}x_{k}^{n-q-\ell}\det\left(
A_{\sim k,\sim q}\right)  =\delta_{\ell,n-q}\left(  -1\right)  ^{n-q}%
\mathfrak{e}_{q-1}\det A.
\]
Hence, Lemma \ref{lem.sol.vander-hook.y-sum.2} is proven in Case 2.

We have now proven Lemma \ref{lem.sol.vander-hook.y-sum.2} in each of the two
Cases 1 and 2. Since these two Cases cover all possibilities, this shows that
Lemma \ref{lem.sol.vander-hook.y-sum.2} always holds.
\end{proof}

\begin{lemma}
\label{lem.sol.vander-hook.lap1-y}Let $n\in\mathbb{N}$. Let $x_{1}%
,x_{2},\ldots,x_{n}$ be $n$ elements of $\mathbb{K}$. Let $t\in\mathbb{K}$.
Let $A$ be the $n\times n$-matrix $\left(  x_{i}^{n-j}\right)  _{1\leq i\leq
n,\ 1\leq j\leq n}$. Let $k\in\left\{  1,2,\ldots,n\right\}  $. Then,%
\begin{align*}
&  V\left(  x_{1},x_{2},\ldots,x_{k-1},x_{k}+t,x_{k+1},x_{k+2},\ldots
,x_{n}\right) \\
&  =\sum_{q=1}^{n}\ \ \sum_{\ell=0}^{n-q}\dbinom{n-q}{\ell}t^{\ell}\left(
-1\right)  ^{k+q}x_{k}^{n-q-\ell}\det\left(  A_{\sim k,\sim q}\right)  .
\end{align*}

\end{lemma}

\begin{proof}
[Proof of Lemma \ref{lem.sol.vander-hook.lap1-y}.]The definition of $A$ yields
$A=\left(  x_{i}^{n-j}\right)  _{1\leq i\leq n,\ 1\leq j\leq n}$.

For every $a\in\mathbb{K}$, $b \in\mathbb{K}$ and $m\in\mathbb{N}$, we have%
\begin{equation}
\left(  a+b\right)  ^{m}=\sum_{\ell=0}^{m}\dbinom{m}{\ell}a^{\ell}b^{m-\ell} .
\label{pf.lem.sol.vander-hook.lap1-y.binom}%
\end{equation}
(Indeed, this is precisely the statement of (\ref{eq.rings.(a+b)**n}), with
the variables $n$ and $k$ renamed as $m$ and $\ell$.) Now, every $q\in\left\{
1,2,\ldots,n\right\}  $ satisfies%
\begin{equation}
\left(  x_{k}+t\right)  ^{n-q}=\sum_{\ell=0}^{n-q}\dbinom{n-q}{\ell}t^{\ell
}x_{k}^{n-q-\ell} \label{pf.lem.sol.vander-hook.lap1-y.binom-used}%
\end{equation}
\footnote{\textit{Proof of (\ref{pf.lem.sol.vander-hook.lap1-y.binom-used}):}
Let $q\in\left\{  1,2,\ldots,n\right\}  $. Thus, $q\leq n$.
\par
Set $m=n-q$. Then, $m=n-q\in\mathbb{N}$ (since $q\leq n$). Applying
(\ref{pf.lem.sol.vander-hook.lap1-y.binom}) to $a=t$ and $b=x_{k}$, we obtain%
\[
\left(  t+x_{k}\right)  ^{m}=\sum_{\ell=0}^{m}\dbinom{m}{\ell}t^{\ell}%
x_{k}^{m-\ell}.
\]
Since $t+x_{k}=x_{k}+t$, this rewrites as
\[
\left(  x_{k}+t\right)  ^{m}=\sum_{\ell=0}^{m}\dbinom{m}{\ell}t^{\ell}%
x_{k}^{m-\ell}.
\]
Since $m=n-q$, this rewrites as
\[
\left(  x_{k}+t\right)  ^{n-q}=\sum_{\ell=0}^{n-q}\dbinom{n-q}{\ell}t^{\ell
}x_{k}^{n-q-\ell}.
\]
This proves (\ref{pf.lem.sol.vander-hook.lap1-y.binom-used}).}.

But%
\begin{align*}
&  V\left(  x_{1},x_{2},\ldots,x_{k-1},x_{k}+t,x_{k+1},x_{k+2},\ldots
,x_{n}\right) \\
&  =\sum_{q=1}^{n}\left(  -1\right)  ^{k+q}\underbrace{\left(  x_{k}+t\right)
^{n-q}}_{\substack{=\sum_{\ell=0}^{n-q}\dbinom{n-q}{\ell}t^{\ell}%
x_{k}^{n-q-\ell}\\\text{(by (\ref{pf.lem.sol.vander-hook.lap1-y.binom-used}%
))}}}\det\left(  A_{\sim k,\sim q}\right)  \ \ \ \ \ \ \ \ \ \ \left(
\text{by Lemma \ref{lem.sol.vander-hook.lap1} \textbf{(b)}}\right) \\
&  =\sum_{q=1}^{n}\left(  -1\right)  ^{k+q}\left(  \sum_{\ell=0}^{n-q}%
\dbinom{n-q}{\ell}t^{\ell}x_{k}^{n-q-\ell}\right)  \det\left(  A_{\sim k,\sim
q}\right) \\
&  =\sum_{q=1}^{n}\ \ \sum_{\ell=0}^{n-q}\underbrace{\left(  -1\right)
^{k+q}\dbinom{n-q}{\ell}t^{\ell}}_{=\dbinom{n-q}{\ell}t^{\ell}\left(
-1\right)  ^{k+q}}x_{k}^{n-q-\ell}\det\left(  A_{\sim k,\sim q}\right) \\
&  =\sum_{q=1}^{n}\ \ \sum_{\ell=0}^{n-q}\dbinom{n-q}{\ell}t^{\ell}\left(
-1\right)  ^{k+q}x_{k}^{n-q-\ell}\det\left(  A_{\sim k,\sim q}\right)  .
\end{align*}
This proves Lemma \ref{lem.sol.vander-hook.lap1-y}.
\end{proof}

\begin{proof}
[Proof of Proposition \ref{prop.sol.vander-hook.variant-1}.]For every
$j\in\mathbb{N}$, define an element $\mathfrak{e}_{j}\in\mathbb{K}$ by
\newline$\mathfrak{e}_{j}=e_{j}\left(  x_{1},x_{2},\ldots,x_{n}\right)  $.

Let $A$ be the $n\times n$-matrix $\left(  x_{i}^{n-j}\right)  _{1\leq i\leq
n,\ 1\leq j\leq n}$. We have%
\begin{align}
&  \sum_{k=1}^{n}y_{k}\underbrace{V\left(  x_{1},x_{2},\ldots,x_{k-1}%
,x_{k}+t,x_{k+1},x_{k+2},\ldots,x_{n}\right)  }_{\substack{=\sum_{q=1}%
^{n}\ \ \sum_{\ell=0}^{n-q}\dbinom{n-q}{\ell}t^{\ell}\left(  -1\right)
^{k+q}x_{k}^{n-q-\ell}\det\left(  A_{\sim k,\sim q}\right)  \\\text{(by Lemma
\ref{lem.sol.vander-hook.lap1-y})}}}\nonumber\\
&  =\sum_{k=1}^{n}\underbrace{y_{k}\sum_{q=1}^{n}\ \ \sum_{\ell=0}%
^{n-q}\dbinom{n-q}{\ell}t^{\ell}\left(  -1\right)  ^{k+q}x_{k}^{n-q-\ell}%
\det\left(  A_{\sim k,\sim q}\right)  }_{=\sum_{q=1}^{n}\ \ \sum_{\ell
=0}^{n-q}y_{k}\dbinom{n-q}{\ell}t^{\ell}\left(  -1\right)  ^{k+q}%
x_{k}^{n-q-\ell}\det\left(  A_{\sim k,\sim q}\right)  }\nonumber\\
&  =\underbrace{\sum_{k=1}^{n}\ \ \sum_{q=1}^{n}\ \ \sum_{\ell=0}^{n-q}%
}_{=\sum_{q=1}^{n}\ \ \sum_{\ell=0}^{n-q}\ \ \sum_{k=1}^{n}}\underbrace{y_{k}%
\dbinom{n-q}{\ell}t^{\ell}\left(  -1\right)  ^{k+q}}_{=\dbinom{n-q}{\ell
}t^{\ell}\left(  -1\right)  ^{k+q}y_{k}}x_{k}^{n-q-\ell}\det\left(  A_{\sim
k,\sim q}\right) \nonumber\\
&  =\sum_{q=1}^{n}\ \ \sum_{\ell=0}^{n-q}\ \ \sum_{k=1}^{n}\dbinom{n-q}{\ell
}t^{\ell}\left(  -1\right)  ^{k+q}y_{k}x_{k}^{n-q-\ell}\det\left(  A_{\sim
k,\sim q}\right)  . \label{pf.prop.sol.vander-hook.variant-1.3}%
\end{align}

Now, let $q\in\left\{  1,2,\ldots,n\right\}  $. Thus, $n-q\in\left\{
0,1,\ldots,n-1\right\}  \subseteq\mathbb{N}$. Hence, $n-q\in\left\{
0,1,\ldots,n-q\right\}  $. Furthermore, (\ref{eq.binom.mm}) (applied to
$m=n-q$) yields $\dbinom{n-q}{n-q}=1$. Now,
\begin{align}
&  \sum_{\ell=0}^{n-q}\ \ \underbrace{\sum_{k=1}^{n}\dbinom{n-q}{\ell}t^{\ell
}\left(  -1\right)  ^{k+q}y_{k}x_{k}^{n-q-\ell}\det\left(  A_{\sim k,\sim
q}\right)  }_{=\dbinom{n-q}{\ell}t^{\ell}\sum_{k=1}^{n}\left(  -1\right)
^{k+q}y_{k}x_{k}^{n-q-\ell}\det\left(  A_{\sim k,\sim q}\right)  }\nonumber\\
&  =\sum_{\ell=0}^{n-q}\dbinom{n-q}{\ell}t^{\ell}\underbrace{\sum_{k=1}%
^{n}\left(  -1\right)  ^{k+q}y_{k}x_{k}^{n-q-\ell}\det\left(  A_{\sim k,\sim
q}\right)  }_{\substack{=\delta_{\ell,n-q}\left(  -1\right)  ^{n-q}%
\mathfrak{e}_{q-1}\det A\\\text{(by Lemma \ref{lem.sol.vander-hook.y-sum.2})}%
}}\nonumber\\
&  =\sum_{\ell=0}^{n-q}\dbinom{n-q}{\ell}t^{\ell}\delta_{\ell,n-q}\left(
-1\right)  ^{n-q}\mathfrak{e}_{q-1}\det A\nonumber\\
&  =\sum_{\ell=0}^{\left(  n-q\right)  -1}\dbinom{n-q}{\ell}t^{\ell
}\underbrace{\delta_{\ell,n-q}}_{\substack{=0\\\text{(because }\ell\neq
n-q\\\text{(since }\ell\leq\left(  n-q\right)  -1<n-q\text{))}}}\left(
-1\right)  ^{n-q}\mathfrak{e}_{q-1}\det A\nonumber\\
&  \ \ \ \ \ \ \ \ \ \ +\underbrace{\dbinom{n-q}{n-q}}_{=1}t^{n-q}%
\underbrace{\delta_{n-q,n-q}}_{\substack{=1\\\text{(since }n-q=n-q\text{)}%
}}\left(  -1\right)  ^{n-q}\mathfrak{e}_{q-1}\det A\nonumber\\
&  \ \ \ \ \ \ \ \ \ \ \left(
\begin{array}
[c]{c}%
\text{here, we have split off the addend for }\ell=n-q\text{ from the sum}\\
\text{(since }n-q\in\left\{  0,1,\ldots,n-q\right\}  \text{)}%
\end{array}
\right) \nonumber\\
&  =\underbrace{\sum_{\ell=0}^{\left(  n-q\right)  -1}\dbinom{n-q}{\ell
}t^{\ell}0\left(  -1\right)  ^{n-q}\mathfrak{e}_{q-1}\det A}_{=0}%
+t^{n-q}\left(  -1\right)  ^{n-q}\mathfrak{e}_{q-1}\det A\nonumber\\
&  =t^{n-q}\left(  -1\right)  ^{n-q}\mathfrak{e}_{q-1}\det A.
\label{pf.prop.sol.vander-hook.variant-1.7}%
\end{align}

Now, forget that we fixed $q$. We thus have proven
(\ref{pf.prop.sol.vander-hook.variant-1.7}) for every $q\in\left\{
1,2,\ldots,n\right\}  $.

Now, (\ref{pf.prop.sol.vander-hook.variant-1.3}) becomes%
\begin{align}
&  \sum_{k=1}^{n}y_{k}V\left(  x_{1},x_{2},\ldots,x_{k-1},x_{k}+t,x_{k+1}%
,x_{k+2},\ldots,x_{n}\right) \nonumber\\
&  =\sum_{q=1}^{n}\ \ \underbrace{\sum_{\ell=0}^{n-q}\ \ \sum_{k=1}^{n}%
\dbinom{n-q}{\ell}t^{\ell}\left(  -1\right)  ^{k+q}y_{k}x_{k}^{n-q-\ell}%
\det\left(  A_{\sim k,\sim q}\right)  }_{\substack{=t^{n-q}\left(  -1\right)
^{n-q}\mathfrak{e}_{q-1}\det A\\\text{(by
(\ref{pf.prop.sol.vander-hook.variant-1.7}))}}}\nonumber\\
&  =\sum_{q=1}^{n}t^{n-q}\left(  -1\right)  ^{n-q}\mathfrak{e}_{q-1}\det A.
\label{pf.prop.sol.vander-hook.variant-1.8}%
\end{align}
But%
\begin{align}
&  \sum_{q=1}^{n}\underbrace{t^{n-q}\left(  -1\right)  ^{n-q}}_{=\left(
t\cdot\left(  -1\right)  \right)  ^{n-q}}\mathfrak{e}_{q-1}\det A\nonumber\\
&  =\sum_{q=1}^{n}\left(  \underbrace{t\cdot\left(  -1\right)  }_{=-t}\right)
^{n-q}\mathfrak{e}_{q-1}\det A=\sum_{q=1}^{n}\left(  -t\right)  ^{n-q}%
\mathfrak{e}_{q-1}\det A\nonumber\\
&  =\sum_{j=0}^{n-1}\underbrace{\left(  -t\right)  ^{n-\left(  n-j\right)  }%
}_{\substack{=\left(  -t\right)  ^{j}\\\text{(since }n-\left(  n-j\right)
=j\text{)}}}\underbrace{\mathfrak{e}_{\left(  n-j\right)  -1}}%
_{\substack{=\mathfrak{e}_{n-1-j}\\\text{(since }\left(  n-j\right)
-1=n-1-j\text{)}}}\det A\nonumber\\
&  \ \ \ \ \ \ \ \ \ \ \left(  \text{here, we have substituted }n-j\text{ for
}q\text{ in the sum}\right) \nonumber\\
&  =\sum_{j=0}^{n-1}\left(  -t\right)  ^{j}\underbrace{\mathfrak{e}_{n-1-j}%
}_{\substack{=e_{n-1-j}\left(  x_{1},x_{2},\ldots,x_{n}\right)  \\\text{(by
the definition of }\mathfrak{e}_{n-1-j}\text{)}}}\det A=\sum_{j=0}%
^{n-1}\left(  -t\right)  ^{j}e_{n-1-j}\left(  x_{1},x_{2},\ldots,x_{n}\right)
\det A\nonumber\\
&  =\left(  \sum_{j=0}^{n-1}\left(  -t\right)  ^{j}e_{n-1-j}\left(
x_{1},x_{2},\ldots,x_{n}\right)  \right)  \det A.
\label{pf.prop.sol.vander-hook.variant-1.9a}%
\end{align}

But the definition of $z_{-t}\left(  x_{1},x_{2},\ldots,x_{n}\right)  $ yields%
\begin{align}
z_{-t}\left(  x_{1},x_{2},\ldots,x_{n}\right)   &  =\sum_{j=0}^{n-1}%
e_{n-1-j}\left(  x_{1},x_{2},\ldots,x_{n}\right)  \left(  -t\right)
^{j}\nonumber\\
&  =\sum_{j=0}^{n-1}\left(  -t\right)  ^{j}e_{n-1-j}\left(  x_{1},x_{2}%
,\ldots,x_{n}\right)  . \label{pf.prop.sol.vander-hook.variant-1.11}%
\end{align}
Thus, (\ref{pf.prop.sol.vander-hook.variant-1.9a}) becomes%
\begin{align}
&  \sum_{q=1}^{n}t^{n-q}\left(  -1\right)  ^{n-q}\mathfrak{e}_{q-1}\det
A\nonumber\\
&  =\underbrace{\left(  \sum_{j=0}^{n-1}\left(  -t\right)  ^{j}e_{n-1-j}%
\left(  x_{1},x_{2},\ldots,x_{n}\right)  \right)  }_{\substack{=z_{-t}\left(
x_{1},x_{2},\ldots,x_{n}\right)  \\\text{(by
(\ref{pf.prop.sol.vander-hook.variant-1.11}))}}}\underbrace{\det
A}_{\substack{=V\left(  x_{1},x_{2},\ldots,x_{n}\right)  \\\text{(by
(\ref{eq.lem.sol.vander-hook.V=x.eq}))}}}\nonumber\\
&  =z_{-t}\left(  x_{1},x_{2},\ldots,x_{n}\right)  \cdot V\left(  x_{1}%
,x_{2},\ldots,x_{n}\right)  . \label{pf.prop.sol.vander-hook.variant-1.13}%
\end{align}
Hence, (\ref{pf.prop.sol.vander-hook.variant-1.8}) becomes%
\begin{align*}
&  \sum_{k=1}^{n}y_{k}V\left(  x_{1},x_{2},\ldots,x_{k-1},x_{k}+t,x_{k+1}%
,x_{k+2},\ldots,x_{n}\right) \\
&  =\sum_{q=1}^{n}t^{n-q}\left(  -1\right)  ^{n-q}\mathfrak{e}_{q-1}\det
A=z_{-t}\left(  x_{1},x_{2},\ldots,x_{n}\right)  \cdot V\left(  x_{1}%
,x_{2},\ldots,x_{n}\right)
\end{align*}
(by (\ref{pf.prop.sol.vander-hook.variant-1.13})). Thus, Proposition
\ref{prop.sol.vander-hook.variant-1} is finally proven.
\end{proof}

Let us finish with one further identity, which is reminiscent of both
Proposition \ref{prop.sol.vander-hook.variant-1} and Proposition
\ref{prop.sol.vander-hook.N}:

\begin{proposition}
\label{prop.sol.vander-hook.variant-2}Let $n\in\mathbb{N}$. Let $x_{1}%
,x_{2},\ldots,x_{n}$ be $n$ elements of $\mathbb{K}$. Let $t\in\mathbb{K}$.
For each $i\in\left\{  1,2,\ldots,n\right\}  $, set $y_{i}=\prod
_{\substack{j\in\left\{  1,2,\ldots,n\right\}  ;\\j\neq i}}x_{j}$. Then,%
\begin{align*}
&  \sum_{k=1}^{n}y_{k}V\left(  x_{1},x_{2},\ldots,x_{k-1},t,x_{k+1}%
,x_{k+2},\ldots,x_{n}\right) \\
&  =z_{-t}\left(  x_{1},x_{2},\ldots,x_{n}\right)  \cdot V\left(  x_{1}%
,x_{2},\ldots,x_{n}\right)  .
\end{align*}

\end{proposition}

The proof proceeds similarly as our proof of Proposition
\ref{prop.sol.vander-hook.N}:

\begin{proof}
[Proof of Proposition \ref{prop.sol.vander-hook.variant-2}.]For every
$j\in\mathbb{N}$, define an element $\mathfrak{e}_{j}\in\mathbb{K}$ by
\newline$\mathfrak{e}_{j}=e_{j}\left(  x_{1},x_{2},\ldots,x_{n}\right)  $.

Let $A$ be the $n\times n$-matrix $\left(  x_{i}^{n-j}\right)  _{1\leq i\leq
n,\ 1\leq j\leq n}$. Now,%
\begin{align*}
&  \sum_{k=1}^{n}y_{k}\underbrace{V\left(  x_{1},x_{2},\ldots,x_{k-1}%
,t,x_{k+1},x_{k+2},\ldots,x_{n}\right)  }_{\substack{=\sum_{q=1}^{n}\left(
-1\right)  ^{k+q}t^{n-q}\det\left(  A_{\sim k,\sim q}\right)  \\\text{(by
Lemma \ref{lem.sol.vander-hook.N.lap1b})}}}\\
&  =\sum_{k=1}^{n}\underbrace{y_{k}\sum_{q=1}^{n}\left(  -1\right)
^{k+q}t^{n-q}\det\left(  A_{\sim k,\sim q}\right)  }_{=\sum_{q=1}^{n}%
y_{k}\left(  -1\right)  ^{k+q}t^{n-q}\det\left(  A_{\sim k,\sim q}\right)  }\\
&  =\underbrace{\sum_{k=1}^{n}\ \ \sum_{q=1}^{n}}_{=\sum_{q=1}^{n}%
\ \ \sum_{k=1}^{n}}\underbrace{y_{k}\left(  -1\right)  ^{k+q}}_{=\left(
-1\right)  ^{k+q}y_{k}}t^{n-q}\det\left(  A_{\sim k,\sim q}\right) \\
&  =\sum_{q=1}^{n}\ \ \underbrace{\sum_{k=1}^{n}\left(  -1\right)  ^{k+q}%
y_{k}t^{n-q}\det\left(  A_{\sim k,\sim q}\right)  }_{=t^{n-q}\sum_{k=1}%
^{n}\left(  -1\right)  ^{k+q}y_{k}\det\left(  A_{\sim k,\sim q}\right)  }\\
&  =\sum_{q=1}^{n}t^{n-q}\underbrace{\sum_{k=1}^{n}\left(  -1\right)
^{k+q}y_{k}\det\left(  A_{\sim k,\sim q}\right)  }_{\substack{=\left(
-1\right)  ^{n-q}\mathfrak{e}_{q-1}\det A\\\text{(by Lemma
\ref{lem.sol.vander-hook.y-sum.1})}}}\\
&  =\sum_{q=1}^{n}t^{n-q}\left(  -1\right)  ^{n-q}\mathfrak{e}_{q-1}\det
A=z_{-t}\left(  x_{1},x_{2},\ldots,x_{n}\right)  \cdot V\left(  x_{1}%
,x_{2},\ldots,x_{n}\right)
\end{align*}
(by (\ref{pf.prop.sol.vander-hook.variant-1.13})). This proves Proposition
\ref{prop.sol.vander-hook.variant-2}.
\end{proof}

\subsection{Solution to Exercise \ref{exe.block2x2.VB+WD}}

\begin{proof}
[Solution to Exercise \ref{exe.block2x2.VB+WD}.]Note that $VB$ and $WD$ are
$m\times m$-matrices. Hence, $VB+WD$ is an $m\times m$-matrix. Also, $\left(
\begin{array}
[c]{cc}%
I_{n} & 0_{n\times m}\\
V & W
\end{array}
\right)  $ and $\left(
\begin{array}
[c]{cc}%
A & B\\
C & D
\end{array}
\right)  $ are $\left(  n+m\right)  \times\left(  n+m\right)  $-matrices.

Exercise \ref{exe.block2x2.mult} (applied to $n$, $m$, $n$, $m$, $n$, $m$,
$I_{n}$, $0_{n\times m}$, $V$, $W$, $A$, $B$, $C$ and $D$ instead of $n$,
$n^{\prime}$, $m$, $m^{\prime}$, $\ell$, $\ell^{\prime}$, $A$, $B$, $C$, $D$,
$A^{\prime}$, $B^{\prime}$, $C^{\prime}$ and $D^{\prime}$) yields%
\begin{align*}
\left(
\begin{array}
[c]{cc}%
I_{n} & 0_{n\times m}\\
V & W
\end{array}
\right)  \left(
\begin{array}
[c]{cc}%
A & B\\
C & D
\end{array}
\right)   &  =\left(
\begin{array}
[c]{cc}%
I_{n}A+0_{n\times m}C & I_{n}B+0_{n\times m}D\\
VA+WC & VB+WD
\end{array}
\right) \\
&  =\left(
\begin{array}
[c]{cc}%
A & B\\
0_{m\times n} & VB+WD
\end{array}
\right)
\end{align*}
(since $I_{n}A+\underbrace{0_{n\times m}C}_{=0_{n\times n}}=I_{n}A=A$,
$I_{n}B+\underbrace{0_{n\times m}D}_{=0_{n\times m}}=I_{n}B=B$ and
$\underbrace{VA}_{=-WC}+WC=-WC+WC=0_{m\times n}$). Taking determinants on both
sides of this equality, we obtain%
\begin{align*}
\det\left(  \left(
\begin{array}
[c]{cc}%
I_{n} & 0_{n\times m}\\
V & W
\end{array}
\right)  \left(
\begin{array}
[c]{cc}%
A & B\\
C & D
\end{array}
\right)  \right)   &  =\det\left(
\begin{array}
[c]{cc}%
A & B\\
0_{m\times n} & VB+WD
\end{array}
\right) \\
&  =\det A\cdot\det\left(  VB+WD\right)
\end{align*}
(by Exercise \ref{exe.block2x2.tridet} (applied to $VB+WD$ instead of $D$)).
Hence,%
\begin{align*}
&  \det A\cdot\det\left(  VB+WD\right) \\
&  =\det\left(  \left(
\begin{array}
[c]{cc}%
I_{n} & 0_{n\times m}\\
V & W
\end{array}
\right)  \left(
\begin{array}
[c]{cc}%
A & B\\
C & D
\end{array}
\right)  \right) \\
&  =\underbrace{\det\left(
\begin{array}
[c]{cc}%
I_{n} & 0_{n\times m}\\
V & W
\end{array}
\right)  }_{\substack{=\det\left(  I_{n}\right)  \cdot\det W\\\text{(by
Exercise \ref{exe.block2x2.tridet.transposed} (applied}\\\text{to }%
I_{n}\text{, }V\text{ and }W\text{ instead of }A\text{, }C\text{ and
}D\text{))}}}\cdot\det\left(
\begin{array}
[c]{cc}%
A & B\\
C & D
\end{array}
\right) \\
&  \ \ \ \ \ \ \ \ \ \ \left(
\begin{array}
[c]{c}%
\text{by Theorem \ref{thm.det(AB)} (applied to }n+m\text{, }\left(
\begin{array}
[c]{cc}%
I_{n} & 0_{n\times m}\\
V & W
\end{array}
\right) \\
\text{and }\left(
\begin{array}
[c]{cc}%
A & B\\
C & D
\end{array}
\right)  \text{ instead of }n\text{, }A\text{ and }B\text{)}%
\end{array}
\right) \\
&  =\underbrace{\det\left(  I_{n}\right)  }_{=1}\cdot\det W\cdot\det\left(
\begin{array}
[c]{cc}%
A & B\\
C & D
\end{array}
\right)  =\det W\cdot\det\left(
\begin{array}
[c]{cc}%
A & B\\
C & D
\end{array}
\right)  .
\end{align*}
This solves Exercise \ref{exe.block2x2.VB+WD}.
\end{proof}

\subsection{Solution to Exercise \ref{exe.block2x2.schur}}

\begin{proof}
[Solution to Exercise \ref{exe.block2x2.schur}.]Recall that $\det\left(
I_{n}\right)  =1$. The same argument (applied to $m$ instead of $n$) yields
$\det\left(  I_{m}\right)  =1$.

The matrix $A$ is invertible; thus, it has an inverse $A^{-1}\in
\mathbb{K}^{n\times n}$. The matrices $I_{m}\in\mathbb{K}^{m\times m}$ and
$-CA^{-1}\in\mathbb{K}^{m\times n}$ satisfy $-\left(  CA^{-1}\right)
A=-I_{m}C$ (since $-\left(  CA^{-1}\right)  A=-C\underbrace{A^{-1}A}_{=I_{n}%
}=-CI_{n}=-\underbrace{C}_{=I_{m}C}=-I_{m}C$). Hence, Exercise
\ref{exe.block2x2.VB+WD} (applied to $W=I_{m}$ and $V=-CA^{-1}$) yields%
\begin{align*}
\det\left(  I_{m}\right)  \cdot\det\left(
\begin{array}
[c]{cc}%
A & B\\
C & D
\end{array}
\right)   &  =\det A\cdot\det\left(  -CA^{-1}B+\underbrace{I_{m}D}_{=D}\right)
\\
&  =\det A\cdot\det\left(  \underbrace{-CA^{-1}B+D}_{=D-CA^{-1}B}\right)
=\det A\cdot\det\left(  D-CA^{-1}B\right)  .
\end{align*}
Hence,%
\[
\det A\cdot\det\left(  D-CA^{-1}B\right)  =\underbrace{\det\left(
I_{m}\right)  }_{=1}\cdot\det\left(
\begin{array}
[c]{cc}%
A & B\\
C & D
\end{array}
\right)  =\det\left(
\begin{array}
[c]{cc}%
A & B\\
C & D
\end{array}
\right)  .
\]
This solves Exercise \ref{exe.block2x2.schur}.
\end{proof}

\subsection{Solution to Exercise \ref{exe.block2x2.jacobi}}

Before we solve Exercise \ref{exe.block2x2.jacobi}, let us show two really
simple facts:

\begin{lemma}
\label{lem.sol.block2x2.jacobi.InIm}Let $n\in\mathbb{N}$ and $m\in\mathbb{N}$.
Then,%
\[
\left(
\begin{array}
[c]{cc}%
I_{n} & 0_{n\times m}\\
0_{m\times n} & I_{m}%
\end{array}
\right)  =I_{n+m}.
\]

\end{lemma}

\begin{vershort}
\begin{proof}
[Proof of Lemma \ref{lem.sol.block2x2.jacobi.InIm}.]Easy, and left to the reader.
\end{proof}
\end{vershort}

\begin{verlong}
\begin{proof}
[Proof of Lemma \ref{lem.sol.block2x2.jacobi.InIm}.]For any two objects $i$
and $j$, define an element $\delta_{i,j}\in\mathbb{K}$ by $\delta_{i,j}=%
\begin{cases}
1, & \text{if }i=j;\\
0, & \text{if }i\neq j
\end{cases}
$.

We have $I_{n}=\left(  \delta_{i,j}\right)  _{1\leq i\leq n,\ 1\leq j\leq n}$
(by the definition of $I_{n}$) and $0_{n\times m}=\left(  0\right)  _{1\leq
i\leq n,\ 1\leq j\leq m}$ (by the definition of $0_{n\times m}$) and
$0_{m\times n}=\left(  0\right)  _{1\leq i\leq m,\ 1\leq j\leq n}$ (by the
definition of $0_{m\times n}$) and $I_{m}=\left(  \delta_{i,j}\right)  _{1\leq
i\leq m,\ 1\leq j\leq m}$ (by the definition of $I_{m}$). Hence, the
definition of $\left(
\begin{array}
[c]{cc}%
I_{n} & 0_{n\times m}\\
0_{m\times n} & I_{m}%
\end{array}
\right)  $ yields%
\begin{align}
&  \left(
\begin{array}
[c]{cc}%
I_{n} & 0_{n\times m}\\
0_{m\times n} & I_{m}%
\end{array}
\right) \nonumber\\
&  =\left(
\begin{cases}
\delta_{i,j}, & \text{if }i\leq n\text{ and }j\leq n;\\
0, & \text{if }i\leq n\text{ and }j>n;\\
0, & \text{if }i>n\text{ and }j\leq n;\\
\delta_{i-n,j-n}, & \text{if }i>n\text{ and }j>n
\end{cases}
\right)  _{1\leq i\leq n+m,\ 1\leq j\leq n+m}.
\label{pf.lem.sol.block2x2.jacobi.InIm.LHS}%
\end{align}

But every $\left(  i,j\right)  \in\left\{  1,2,\ldots,n+m\right\}  ^{2}$
satisfies%
\begin{equation}%
\begin{cases}
\delta_{i,j}, & \text{if }i\leq n\text{ and }j\leq n;\\
0, & \text{if }i\leq n\text{ and }j>n;\\
0, & \text{if }i>n\text{ and }j\leq n;\\
\delta_{i-n,j-n}, & \text{if }i>n\text{ and }j>n
\end{cases}
=\delta_{i,j} \label{pf.lem.sol.block2x2.jacobi.InIm.deldel}%
\end{equation}
\footnote{\textit{Proof of (\ref{pf.lem.sol.block2x2.jacobi.InIm.deldel}):}
Let $\left(  i,j\right)  \in\left\{  1,2,\ldots,n+m\right\}  ^{2}$. We want to
prove (\ref{pf.lem.sol.block2x2.jacobi.InIm.deldel}).
\par
We have $\left(  i,j\right)  \in\left\{  1,2,\ldots,n+m\right\}  ^{2}$. In
other words, $i\in\left\{  1,2,\ldots,n+m\right\}  $ and $j\in\left\{
1,2,\ldots,n+m\right\}  $. We are in one of the following two cases:
\par
\textit{Case 1:} We have $i\leq n$.
\par
\textit{Case 2:} We have $i>n$.
\par
Let us first consider Case 1. In this case, we have $i\leq n$. We are in one
of the following two subcases:
\par
\textit{Subcase 1.1:} We have $j\leq n$.
\par
\textit{Subcase 1.2:} We have $j>n$.
\par
Let us first consider Subcase 1.1. In this case, we have $j\leq n$. Hence,
\[%
\begin{cases}
\delta_{i,j}, & \text{if }i\leq n\text{ and }j\leq n;\\
0, & \text{if }i\leq n\text{ and }j>n;\\
0, & \text{if }i>n\text{ and }j\leq n;\\
\delta_{i-n,j-n}, & \text{if }i>n\text{ and }j>n
\end{cases}
=\delta_{i,j}\ \ \ \ \ \ \ \ \ \ \left(  \text{since }i\leq n\text{ and }j\leq
n\right)  .
\]
Thus, (\ref{pf.lem.sol.block2x2.jacobi.InIm.deldel}) is proven in Subcase 1.1.
\par
Let us next consider Subcase 1.2. In this case, we have $j>n$. Hence, $n<j$,
so that $i\leq n<j$. Thus, $i\neq j$, so that $\delta_{i,j}=0$. But%
\begin{align*}%
\begin{cases}
\delta_{i,j}, & \text{if }i\leq n\text{ and }j\leq n;\\
0, & \text{if }i\leq n\text{ and }j>n;\\
0, & \text{if }i>n\text{ and }j\leq n;\\
\delta_{i-n,j-n}, & \text{if }i>n\text{ and }j>n
\end{cases}
&  =0\ \ \ \ \ \ \ \ \ \ \left(  \text{since }i\leq n\text{ and }j>n\right) \\
&  =\delta_{i,j}\ \ \ \ \ \ \ \ \ \ \left(  \text{since }\delta_{i,j}%
=0\right)  .
\end{align*}
Thus, (\ref{pf.lem.sol.block2x2.jacobi.InIm.deldel}) is proven in Subcase 1.2.
\par
We now have proven (\ref{pf.lem.sol.block2x2.jacobi.InIm.deldel}) in each of
the two Subcases 1.1 and 1.2. Since these two Subcases cover the whole Case 1,
we can thus conclude that (\ref{pf.lem.sol.block2x2.jacobi.InIm.deldel}) holds
in Case 1.
\par
Let us now consider Case 2. In this case, we have $i>n$. We are in one of the
following two subcases:
\par
\textit{Subcase 2.1:} We have $j\leq n$.
\par
\textit{Subcase 2.2:} We have $j>n$.
\par
Let us first consider Subcase 2.1. In this case, we have $j\leq n$. Hence,
$n\geq j$, and thus $i>n\geq j$. Hence, $i\neq j$, and thus $\delta_{i,j}=0$.
But%
\begin{align*}%
\begin{cases}
\delta_{i,j}, & \text{if }i\leq n\text{ and }j\leq n;\\
0, & \text{if }i\leq n\text{ and }j>n;\\
0, & \text{if }i>n\text{ and }j\leq n;\\
\delta_{i-n,j-n}, & \text{if }i>n\text{ and }j>n
\end{cases}
&  =0\ \ \ \ \ \ \ \ \ \ \left(  \text{since }i>n\text{ and }j\leq n\right) \\
&  =\delta_{i,j}\ \ \ \ \ \ \ \ \ \ \left(  \text{since }\delta_{i,j}%
=0\right)  .
\end{align*}
Thus, (\ref{pf.lem.sol.block2x2.jacobi.InIm.deldel}) is proven in Subcase 2.1.
\par
Let us next consider Subcase 2.2. In this case, we have $j>n$. But the
definition of $\delta_{i,j}$ yields $\delta_{i,j}=%
\begin{cases}
1, & \text{if }i=j;\\
0, & \text{if }i\neq j
\end{cases}
$. Now,%
\begin{align*}
&
\begin{cases}
\delta_{i,j}, & \text{if }i\leq n\text{ and }j\leq n;\\
0, & \text{if }i\leq n\text{ and }j>n;\\
0, & \text{if }i>n\text{ and }j\leq n;\\
\delta_{i-n,j-n}, & \text{if }i>n\text{ and }j>n
\end{cases}
\\
&  =\delta_{i-n,j-n}\ \ \ \ \ \ \ \ \ \ \left(  \text{since }i>n\text{ and
}j>n\right) \\
&  =%
\begin{cases}
1, & \text{if }i-n=j-n;\\
0, & \text{if }i-n\neq j-n
\end{cases}
\\
&  =%
\begin{cases}
1, & \text{if }i=j;\\
0, & \text{if }i\neq j
\end{cases}
\\
&  \ \ \ \ \ \ \ \ \ \ \left(
\begin{array}
[c]{c}%
\text{since the statement }i-n=j-n\text{ is equivalent to }i=j\text{,}\\
\text{and since the statement }i-n\neq j-n\text{ is equivalent to }i\neq j
\end{array}
\right) \\
&  =\delta_{i,j}.
\end{align*}
Thus, (\ref{pf.lem.sol.block2x2.jacobi.InIm.deldel}) is proven in Subcase 2.2.
\par
We now have proven (\ref{pf.lem.sol.block2x2.jacobi.InIm.deldel}) in each of
the two Subcases 2.1 and 2.2. Since these two Subcases cover the whole Case 2,
we can thus conclude that (\ref{pf.lem.sol.block2x2.jacobi.InIm.deldel}) holds
in Case 2.
\par
We now have proven (\ref{pf.lem.sol.block2x2.jacobi.InIm.deldel}) in each of
the two Cases 1 and 2. Since these two Cases cover all possibilities, we can
thus conclude that (\ref{pf.lem.sol.block2x2.jacobi.InIm.deldel}) always
holds. Qed.}.

Now, (\ref{pf.lem.sol.block2x2.jacobi.InIm.LHS}) becomes%
\begin{align*}
\left(
\begin{array}
[c]{cc}%
I_{n} & 0_{n\times m}\\
0_{m\times n} & I_{m}%
\end{array}
\right)   &  =\left(  \underbrace{%
\begin{cases}
\delta_{i,j}, & \text{if }i\leq n\text{ and }j\leq n;\\
0, & \text{if }i\leq n\text{ and }j>n;\\
0, & \text{if }i>n\text{ and }j\leq n;\\
\delta_{i-n,j-n}, & \text{if }i>n\text{ and }j>n
\end{cases}
}_{\substack{=\delta_{i,j}\\\text{(by
(\ref{pf.lem.sol.block2x2.jacobi.InIm.deldel}))}}}\right)  _{1\leq i\leq
n+m,\ 1\leq j\leq n+m}\\
&  =\left(  \delta_{i,j}\right)  _{1\leq i\leq n+m,\ 1\leq j\leq n+m}.
\end{align*}
Comparing this with%
\[
I_{n+m}=\left(  \delta_{i,j}\right)  _{1\leq i\leq n+m,\ 1\leq j\leq
n+m}\ \ \ \ \ \ \ \ \ \ \left(  \text{by the definition of }I_{n+m}\right)  ,
\]
we obtain $\left(
\begin{array}
[c]{cc}%
I_{n} & 0_{n\times m}\\
0_{m\times n} & I_{m}%
\end{array}
\right)  =I_{n+m}$. This proves Lemma \ref{lem.sol.block2x2.jacobi.InIm}.
\end{proof}
\end{verlong}

\begin{lemma}
\label{lem.sol.block2x2.jacobi.equal}Let $n$, $n^{\prime}$, $m$ and
$m^{\prime}$ be four nonnegative integers.

Let $A\in\mathbb{K}^{n\times m}$ and $A^{\prime}\in\mathbb{K}^{n\times m}$.
Let $B\in\mathbb{K}^{n\times m^{\prime}}$ and $B^{\prime}\in\mathbb{K}%
^{n\times m^{\prime}}$. Let $C\in\mathbb{K}^{n^{\prime}\times m}$ and
$C^{\prime}\in\mathbb{K}^{n^{\prime}\times m}$. Let $D\in\mathbb{K}%
^{n^{\prime}\times m^{\prime}}$ and $D^{\prime}\in\mathbb{K}^{n^{\prime}\times
m^{\prime}}$. Assume that $\left(
\begin{array}
[c]{cc}%
A & B\\
C & D
\end{array}
\right)  =\left(
\begin{array}
[c]{cc}%
A^{\prime} & B^{\prime}\\
C^{\prime} & D^{\prime}%
\end{array}
\right)  $. Then, $A=A^{\prime}$, $B=B^{\prime}$, $C=C^{\prime}$ and
$D=D^{\prime}$.
\end{lemma}

\begin{vershort}
\begin{proof}
[Proof of Lemma \ref{lem.sol.block2x2.jacobi.equal}.]Easy, and left to the reader.
\end{proof}
\end{vershort}

\begin{verlong}
\begin{proof}
[Proof of Lemma \ref{lem.sol.block2x2.jacobi.equal}.]Write the matrix
$A\in\mathbb{K}^{n\times m}$ in the form $A=\left(  a_{i,j}\right)  _{1\leq
i\leq n,\ 1\leq j\leq m}$.

Write the matrix $A^{\prime}\in\mathbb{K}^{n\times m}$ in the form $A^{\prime
}=\left(  a_{i,j}^{\prime}\right)  _{1\leq i\leq n,\ 1\leq j\leq m}$.

Write the matrix $B\in\mathbb{K}^{n\times m^{\prime}}$ in the form $B=\left(
b_{i,j}\right)  _{1\leq i\leq n,\ 1\leq j\leq m^{\prime}}$.

Write the matrix $B^{\prime}\in\mathbb{K}^{n\times m^{\prime}}$ in the form
$B^{\prime}=\left(  b_{i,j}^{\prime}\right)  _{1\leq i\leq n,\ 1\leq j\leq
m^{\prime}}$.

Write the matrix $C\in\mathbb{K}^{n^{\prime}\times m}$ in the form $C=\left(
c_{i,j}\right)  _{1\leq i\leq n^{\prime},\ 1\leq j\leq m}$.

Write the matrix $C^{\prime}\in\mathbb{K}^{n^{\prime}\times m}$ in the form
$C^{\prime}=\left(  c_{i,j}^{\prime}\right)  _{1\leq i\leq n^{\prime},\ 1\leq
j\leq m}$.

Write the matrix $D\in\mathbb{K}^{n^{\prime}\times m^{\prime}}$ in the form
$D=\left(  d_{i,j}\right)  _{1\leq i\leq n^{\prime},\ 1\leq j\leq m^{\prime}}$.

Write the matrix $D^{\prime}\in\mathbb{K}^{n^{\prime}\times m^{\prime}}$ in
the form $D^{\prime}=\left(  d_{i,j}^{\prime}\right)  _{1\leq i\leq n^{\prime
},\ 1\leq j\leq m^{\prime}}$.

We have $A=\left(  a_{i,j}\right)  _{1\leq i\leq n,\ 1\leq j\leq m}$,
$B=\left(  b_{i,j}\right)  _{1\leq i\leq n,\ 1\leq j\leq m^{\prime}}$,
$C=\left(  c_{i,j}\right)  _{1\leq i\leq n^{\prime},\ 1\leq j\leq m}$ and
$D=\left(  d_{i,j}\right)  _{1\leq i\leq n^{\prime},\ 1\leq j\leq m^{\prime}}%
$. Hence, the definition of $\left(
\begin{array}
[c]{cc}%
A & B\\
C & D
\end{array}
\right)  $ yields%
\begin{align}
&  \left(
\begin{array}
[c]{cc}%
A & B\\
C & D
\end{array}
\right) \nonumber\\
&  =\left(
\begin{cases}
a_{i,j}, & \text{if }i\leq n\text{ and }j\leq m;\\
b_{i,j-m}, & \text{if }i\leq n\text{ and }j>m;\\
c_{i-n,j}, & \text{if }i>n\text{ and }j\leq m;\\
d_{i-n,j-m}, & \text{if }i>n\text{ and }j>m
\end{cases}
\right)  _{1\leq i\leq n+n^{\prime},\ 1\leq j\leq m+m^{\prime}}.
\label{pf.lem.sol.block2x2.jacobi.equal.1}%
\end{align}

We have $A^{\prime}=\left(  a_{i,j}^{\prime}\right)  _{1\leq i\leq n,\ 1\leq
j\leq m}$, $B^{\prime}=\left(  b_{i,j}^{\prime}\right)  _{1\leq i\leq
n,\ 1\leq j\leq m^{\prime}}$, $C^{\prime}=\left(  c_{i,j}^{\prime}\right)
_{1\leq i\leq n^{\prime},\ 1\leq j\leq m}$ and $D^{\prime}=\left(
d_{i,j}^{\prime}\right)  _{1\leq i\leq n^{\prime},\ 1\leq j\leq m^{\prime}}$.
Hence, the definition of $\left(
\begin{array}
[c]{cc}%
A^{\prime} & B^{\prime}\\
C^{\prime} & D^{\prime}%
\end{array}
\right)  $ yields%
\begin{align}
&  \left(
\begin{array}
[c]{cc}%
A^{\prime} & B^{\prime}\\
C^{\prime} & D^{\prime}%
\end{array}
\right) \nonumber\\
&  =\left(
\begin{cases}
a_{i,j}^{\prime}, & \text{if }i\leq n\text{ and }j\leq m;\\
b_{i,j-m}^{\prime}, & \text{if }i\leq n\text{ and }j>m;\\
c_{i-n,j}^{\prime}, & \text{if }i>n\text{ and }j\leq m;\\
d_{i-n,j-m}^{\prime}, & \text{if }i>n\text{ and }j>m
\end{cases}
\right)  _{1\leq i\leq n+n^{\prime},\ 1\leq j\leq m+m^{\prime}}.
\label{pf.lem.sol.block2x2.jacobi.equal.2}%
\end{align}
Now, from (\ref{pf.lem.sol.block2x2.jacobi.equal.1}), we obtain%
\begin{align*}
&  \left(
\begin{cases}
a_{i,j}, & \text{if }i\leq n\text{ and }j\leq m;\\
b_{i,j-m}, & \text{if }i\leq n\text{ and }j>m;\\
c_{i-n,j}, & \text{if }i>n\text{ and }j\leq m;\\
d_{i-n,j-m}, & \text{if }i>n\text{ and }j>m
\end{cases}
\right)  _{1\leq i\leq n+n^{\prime},\ 1\leq j\leq m+m^{\prime}}\\
&  =\left(
\begin{array}
[c]{cc}%
A & B\\
C & D
\end{array}
\right)  =\left(
\begin{array}
[c]{cc}%
A^{\prime} & B^{\prime}\\
C^{\prime} & D^{\prime}%
\end{array}
\right) \\
&  =\left(
\begin{cases}
a_{i,j}^{\prime}, & \text{if }i\leq n\text{ and }j\leq m;\\
b_{i,j-m}^{\prime}, & \text{if }i\leq n\text{ and }j>m;\\
c_{i-n,j}^{\prime}, & \text{if }i>n\text{ and }j\leq m;\\
d_{i-n,j-m}^{\prime}, & \text{if }i>n\text{ and }j>m
\end{cases}
\right)  _{1\leq i\leq n+n^{\prime},\ 1\leq j\leq m+m^{\prime}}%
\end{align*}
(by (\ref{pf.lem.sol.block2x2.jacobi.equal.2})). In other words,%
\begin{align}
&
\begin{cases}
a_{i,j}, & \text{if }i\leq n\text{ and }j\leq m;\\
b_{i,j-m}, & \text{if }i\leq n\text{ and }j>m;\\
c_{i-n,j}, & \text{if }i>n\text{ and }j\leq m;\\
d_{i-n,j-m}, & \text{if }i>n\text{ and }j>m
\end{cases}
\nonumber\\
&  =%
\begin{cases}
a_{i,j}^{\prime}, & \text{if }i\leq n\text{ and }j\leq m;\\
b_{i,j-m}^{\prime}, & \text{if }i\leq n\text{ and }j>m;\\
c_{i-n,j}^{\prime}, & \text{if }i>n\text{ and }j\leq m;\\
d_{i-n,j-m}^{\prime}, & \text{if }i>n\text{ and }j>m
\end{cases}
\label{pf.lem.sol.block2x2.jacobi.equal.3}%
\end{align}
for every $\left(  i,j\right)  \in\left\{  1,2,\ldots,n+n^{\prime}\right\}
\times\left\{  1,2,\ldots,m+m^{\prime}\right\}  $.

Every $\left(  i,j\right)  \in\left\{  1,2,\ldots,n\right\}  \times\left\{
1,2,\ldots,m\right\}  $ satisfies $a_{i,j}=a_{i,j}^{\prime}$%
\ \ \ \ \footnote{\textit{Proof.} Let $\left(  i,j\right)  \in\left\{
1,2,\ldots,n\right\}  \times\left\{  1,2,\ldots,m\right\}  $. Thus,
$i\in\left\{  1,2,\ldots,n\right\}  $ and $j\in\left\{  1,2,\ldots,m\right\}
$.
\par
We have $i\in\left\{  1,2,\ldots,n\right\}  \subseteq\left\{  1,2,\ldots
,n+n^{\prime}\right\}  $ (since $n\leq n+n^{\prime}$) and $i\leq n$ (since
$i\in\left\{  1,2,\ldots,n\right\}  $). We have $j\in\left\{  1,2,\ldots
,m\right\}  \subseteq\left\{  1,2,\ldots,m+m^{\prime}\right\}  $ (since $m\leq
m+m^{\prime}$) and $j\leq m$ (since $j\in\left\{  1,2,\ldots,m\right\}  $).
\par
From $i\in\left\{  1,2,\ldots,n+n^{\prime}\right\}  $ and $j\in\left\{
1,2,\ldots,m+m^{\prime}\right\}  $, we obtain $\left(  i,j\right)  \in\left\{
1,2,\ldots,n+n^{\prime}\right\}  \times\left\{  1,2,\ldots,m+m^{\prime
}\right\}  $. Now, $%
\begin{cases}
a_{i,j}, & \text{if }i\leq n\text{ and }j\leq m;\\
b_{i,j-m}, & \text{if }i\leq n\text{ and }j>m;\\
c_{i-n,j}, & \text{if }i>n\text{ and }j\leq m;\\
d_{i-n,j-m}, & \text{if }i>n\text{ and }j>m
\end{cases}
=a_{i,j}$ (since $i\leq n$ and $j\leq m$). Hence,%
\begin{align*}
a_{i,j}  &  =%
\begin{cases}
a_{i,j}, & \text{if }i\leq n\text{ and }j\leq m;\\
b_{i,j-m}, & \text{if }i\leq n\text{ and }j>m;\\
c_{i-n,j}, & \text{if }i>n\text{ and }j\leq m;\\
d_{i-n,j-m}, & \text{if }i>n\text{ and }j>m
\end{cases}
\\
&  =%
\begin{cases}
a_{i,j}^{\prime}, & \text{if }i\leq n\text{ and }j\leq m;\\
b_{i,j-m}^{\prime}, & \text{if }i\leq n\text{ and }j>m;\\
c_{i-n,j}^{\prime}, & \text{if }i>n\text{ and }j\leq m;\\
d_{i-n,j-m}^{\prime}, & \text{if }i>n\text{ and }j>m
\end{cases}
\ \ \ \ \ \ \ \ \ \ \left(  \text{by (\ref{pf.lem.sol.block2x2.jacobi.equal.3}%
)}\right) \\
&  =a_{i,j}^{\prime}\ \ \ \ \ \ \ \ \ \ \left(  \text{since }i\leq n\text{ and
}j\leq m\right)  .
\end{align*}
Qed.}. In other words, $\left(  a_{i,j}\right)  _{1\leq i\leq n,\ 1\leq j\leq
m}=\left(  a_{i,j}^{\prime}\right)  _{1\leq i\leq n,\ 1\leq j\leq m}$. Now,%
\[
A=\left(  a_{i,j}\right)  _{1\leq i\leq n,\ 1\leq j\leq m}=\left(
a_{i,j}^{\prime}\right)  _{1\leq i\leq n,\ 1\leq j\leq m}=A^{\prime}.
\]

Every $\left(  i,j\right)  \in\left\{  1,2,\ldots,n\right\}  \times\left\{
1,2,\ldots,m^{\prime}\right\}  $ satisfies $b_{i,j}=b_{i,j}^{\prime}%
$\ \ \ \ \footnote{\textit{Proof.} Let $\left(  i,j\right)  \in\left\{
1,2,\ldots,n\right\}  \times\left\{  1,2,\ldots,m^{\prime}\right\}  $. Thus,
$i\in\left\{  1,2,\ldots,n\right\}  $ and $j\in\left\{  1,2,\ldots,m^{\prime
}\right\}  $. We have $i\in\left\{  1,2,\ldots,n\right\}  \subseteq\left\{
1,2,\ldots,n+n^{\prime}\right\}  $ (since $n\leq n+n^{\prime}$) and $i\leq n$
(since $i\in\left\{  1,2,\ldots,n\right\}  $). We have $j\in\left\{
1,2,\ldots,m^{\prime}\right\}  $ and thus $j+m\in\left\{  m+1,m+2,\ldots
,m+m^{\prime}\right\}  \subseteq\left\{  1,2,\ldots,m+m^{\prime}\right\}  $
(since $\underbrace{m}_{\geq0}+1\geq1$). Also, $\underbrace{j}_{>0}+m>m$.
\par
From $i\in\left\{  1,2,\ldots,n+n^{\prime}\right\}  $ and $j+m\in\left\{
1,2,\ldots,m+m^{\prime}\right\}  $, we obtain $\left(  i,j+m\right)
\in\left\{  1,2,\ldots,n+n^{\prime}\right\}  \times\left\{  1,2,\ldots
,m+m^{\prime}\right\}  $. Now, $%
\begin{cases}
a_{i,j}, & \text{if }i\leq n\text{ and }j+m\leq m;\\
b_{i,\left(  j+m\right)  -m}, & \text{if }i\leq n\text{ and }j+m>m;\\
c_{i-n,j+m}, & \text{if }i>n\text{ and }j+m\leq m;\\
d_{i-n,\left(  j+m\right)  -m}, & \text{if }i>n\text{ and }j+m>m
\end{cases}
=b_{i,\left(  j+m\right)  -m}$ (since $i\leq n$ and $j+m>m$). Hence,%
\begin{align*}
b_{i,\left(  j+m\right)  -m}  &  =%
\begin{cases}
a_{i,j}, & \text{if }i\leq n\text{ and }j+m\leq m;\\
b_{i,\left(  j+m\right)  -m}, & \text{if }i\leq n\text{ and }j+m>m;\\
c_{i-n,j+m}, & \text{if }i>n\text{ and }j+m\leq m;\\
d_{i-n,\left(  j+m\right)  -m}, & \text{if }i>n\text{ and }j+m>m
\end{cases}
\\
&  =%
\begin{cases}
a_{i,j}^{\prime}, & \text{if }i\leq n\text{ and }j+m\leq m;\\
b_{i,\left(  j+m\right)  -m}^{\prime}, & \text{if }i\leq n\text{ and }j+m>m;\\
c_{i-n,j+m}^{\prime}, & \text{if }i>n\text{ and }j+m\leq m;\\
d_{i-n,\left(  j+m\right)  -m}^{\prime}, & \text{if }i>n\text{ and }j+m>m
\end{cases}
\ \ \ \ \ \ \ \ \ \ \left(
\begin{array}
[c]{c}%
\text{by (\ref{pf.lem.sol.block2x2.jacobi.equal.3}) (applied}\\
\text{to }\left(  i,j+m\right)  \text{ instead of }\left(  i,j\right)
\text{)}%
\end{array}
\right) \\
&  =b_{i,\left(  j+m\right)  -m}^{\prime}\ \ \ \ \ \ \ \ \ \ \left(
\text{since }i\leq n\text{ and }j+m>m\right)  .
\end{align*}
This rewrites as $b_{i,j}=b_{i,j}^{\prime}$ (since $\left(  j+m\right)
-m=j$). Qed.}. In other words, $\left(  b_{i,j}\right)  _{1\leq i\leq
n,\ 1\leq j\leq m^{\prime}}=\left(  b_{i,j}^{\prime}\right)  _{1\leq i\leq
n,\ 1\leq j\leq m^{\prime}}$. Now,%
\[
B=\left(  b_{i,j}\right)  _{1\leq i\leq n,\ 1\leq j\leq m^{\prime}}=\left(
b_{i,j}^{\prime}\right)  _{1\leq i\leq n,\ 1\leq j\leq m^{\prime}}=B^{\prime
}.
\]

Every $\left(  i,j\right)  \in\left\{  1,2,\ldots,n^{\prime}\right\}
\times\left\{  1,2,\ldots,m\right\}  $ satisfies $c_{i,j}=c_{i,j}^{\prime}%
$\ \ \ \ \footnote{\textit{Proof.} Let $\left(  i,j\right)  \in\left\{
1,2,\ldots,n^{\prime}\right\}  \times\left\{  1,2,\ldots,m\right\}  $. Thus,
$i\in\left\{  1,2,\ldots,n^{\prime}\right\}  $ and $j\in\left\{
1,2,\ldots,m\right\}  $. We have $i\in\left\{  1,2,\ldots,n^{\prime}\right\}
$ and thus $i+n\in\left\{  n+1,n+2,\ldots,n+n^{\prime}\right\}  \subseteq
\left\{  1,2,\ldots,n+n^{\prime}\right\}  $ (since $\underbrace{n}_{\geq
0}+1\geq n$). Also, $\underbrace{i}_{>0}+n>n$. We have $j\in\left\{
1,2,\ldots,m\right\}  \subseteq\left\{  1,2,\ldots,m+m^{\prime}\right\}  $
(since $m\leq m+m^{\prime}$) and $j\leq m$ (since $j\in\left\{  1,2,\ldots
,m\right\}  $).
\par
From $i+n\in\left\{  1,2,\ldots,n+n^{\prime}\right\}  $ and $j\in\left\{
1,2,\ldots,m+m^{\prime}\right\}  $, we obtain $\left(  i+n,j\right)
\in\left\{  1,2,\ldots,n+n^{\prime}\right\}  \times\left\{  1,2,\ldots
,m+m^{\prime}\right\}  $. Now, $%
\begin{cases}
a_{i+n,j}, & \text{if }i+n\leq n\text{ and }j\leq m;\\
b_{i+n,j-m}, & \text{if }i+n\leq n\text{ and }j>m;\\
c_{\left(  i+n\right)  -n,j}, & \text{if }i+n>n\text{ and }j\leq m;\\
d_{\left(  i+n\right)  -n,j-m}, & \text{if }i+n>n\text{ and }j>m
\end{cases}
=c_{\left(  i+n\right)  -n,j}$ (since $i+n>n$ and $j\leq m$). Hence,%
\begin{align*}
c_{\left(  i+n\right)  -n,j}  &  =%
\begin{cases}
a_{i+n,j}, & \text{if }i+n\leq n\text{ and }j\leq m;\\
b_{i+n,j-m}, & \text{if }i+n\leq n\text{ and }j>m;\\
c_{\left(  i+n\right)  -n,j}, & \text{if }i+n>n\text{ and }j\leq m;\\
d_{\left(  i+n\right)  -n,j-m}, & \text{if }i+n>n\text{ and }j>m
\end{cases}
\\
&  =%
\begin{cases}
a_{i+n,j}^{\prime}, & \text{if }i+n\leq n\text{ and }j\leq m;\\
b_{i+n,j-m}^{\prime}, & \text{if }i+n\leq n\text{ and }j>m;\\
c_{\left(  i+n\right)  -n,j}^{\prime}, & \text{if }i+n>n\text{ and }j\leq m;\\
d_{\left(  i+n\right)  -n,j-m}^{\prime}, & \text{if }i+n>n\text{ and }j>m
\end{cases}
\ \ \ \ \ \ \ \ \ \ \left(
\begin{array}
[c]{c}%
\text{by (\ref{pf.lem.sol.block2x2.jacobi.equal.3}) (applied}\\
\text{to }\left(  i+n,j\right)  \text{ instead of }\left(  i,j\right)
\text{)}%
\end{array}
\right) \\
&  =c_{\left(  i+n\right)  -n,j}^{\prime}\ \ \ \ \ \ \ \ \ \ \left(
\text{since }i+n>n\text{ and }j\leq m\right)  .
\end{align*}
This rewrites as $c_{i,j}=c_{i,j}^{\prime}$ (since $\left(  i+n\right)
-n=i$). Qed.}. In other words, $\left(  c_{i,j}\right)  _{1\leq i\leq
n^{\prime},\ 1\leq j\leq m}=\left(  c_{i,j}^{\prime}\right)  _{1\leq i\leq
n^{\prime},\ 1\leq j\leq m}$. Now,%
\[
C=\left(  c_{i,j}\right)  _{1\leq i\leq n^{\prime},\ 1\leq j\leq m}=\left(
c_{i,j}^{\prime}\right)  _{1\leq i\leq n^{\prime},\ 1\leq j\leq m}=C^{\prime
}.
\]

Every $\left(  i,j\right)  \in\left\{  1,2,\ldots,n^{\prime}\right\}
\times\left\{  1,2,\ldots,m^{\prime}\right\}  $ satisfies $d_{i,j}%
=d_{i,j}^{\prime}$\ \ \ \ \footnote{\textit{Proof.} Let $\left(  i,j\right)
\in\left\{  1,2,\ldots,n^{\prime}\right\}  \times\left\{  1,2,\ldots
,m^{\prime}\right\}  $. Thus, $i\in\left\{  1,2,\ldots,n^{\prime}\right\}  $
and $j\in\left\{  1,2,\ldots,m^{\prime}\right\}  $. We have $i\in\left\{
1,2,\ldots,n^{\prime}\right\}  $ and thus $i+n\in\left\{  n+1,n+2,\ldots
,n+n^{\prime}\right\}  \subseteq\left\{  1,2,\ldots,n+n^{\prime}\right\}  $
(since $\underbrace{n}_{\geq0}+1\geq n$). Also, $\underbrace{i}_{>0}+n>n$. We
have $j\in\left\{  1,2,\ldots,m^{\prime}\right\}  $ and thus $j+m\in\left\{
m+1,m+2,\ldots,m+m^{\prime}\right\}  \subseteq\left\{  1,2,\ldots,m+m^{\prime
}\right\}  $ (since $\underbrace{m}_{\geq0}+1\geq m$). Also, $\underbrace{j}%
_{>0}+m>m$.
\par
From $i+n\in\left\{  1,2,\ldots,n+n^{\prime}\right\}  $ and $j+m\in\left\{
1,2,\ldots,m+m^{\prime}\right\}  $, we obtain $\left(  i+n,j+m\right)
\in\left\{  1,2,\ldots,n+n^{\prime}\right\}  \times\left\{  1,2,\ldots
,m+m^{\prime}\right\}  $. Now, $%
\begin{cases}
a_{i+n,j+m}, & \text{if }i+n\leq n\text{ and }j+m\leq m;\\
b_{i+n,\left(  j+m\right)  -m}, & \text{if }i+n\leq n\text{ and }j+m>m;\\
c_{\left(  i+n\right)  -n,j+m}, & \text{if }i+n>n\text{ and }j+m\leq m;\\
d_{\left(  i+n\right)  -n,\left(  j+m\right)  -m}, & \text{if }i+n>n\text{ and
}j+m>m
\end{cases}
=d_{\left(  i+n\right)  -n,\left(  j+m\right)  -m}$ (since $i+n>n$ and
$j+m>m$). Hence,%
\begin{align*}
&  d_{\left(  i+n\right)  -n,\left(  j+m\right)  -m}\\
&  =%
\begin{cases}
a_{i+n,j+m}, & \text{if }i+n\leq n\text{ and }j+m\leq m;\\
b_{i+n,\left(  j+m\right)  -m}, & \text{if }i+n\leq n\text{ and }j+m>m;\\
c_{\left(  i+n\right)  -n,j+m}, & \text{if }i+n>n\text{ and }j+m\leq m;\\
d_{\left(  i+n\right)  -n,\left(  j+m\right)  -m}, & \text{if }i+n>n\text{ and
}j+m>m
\end{cases}
\\
&  =%
\begin{cases}
a_{i+n,j+m}^{\prime}, & \text{if }i+n\leq n\text{ and }j+m\leq m;\\
b_{i+n,\left(  j+m\right)  -m}^{\prime}, & \text{if }i+n\leq n\text{ and
}j+m>m;\\
c_{\left(  i+n\right)  -n,j+m}^{\prime}, & \text{if }i+n>n\text{ and }j+m\leq
m;\\
d_{\left(  i+n\right)  -n,\left(  j+m\right)  -m}^{\prime}, & \text{if
}i+n>n\text{ and }j+m>m
\end{cases}
\ \ \ \ \ \ \ \ \ \ \left(
\begin{array}
[c]{c}%
\text{by (\ref{pf.lem.sol.block2x2.jacobi.equal.3}) (applied}\\
\text{to }\left(  i+n,j+m\right)  \text{ instead of }\left(  i,j\right)
\text{)}%
\end{array}
\right) \\
&  =d_{\left(  i+n\right)  -n,\left(  j+m\right)  -m}^{\prime}%
\ \ \ \ \ \ \ \ \ \ \left(  \text{since }i+n>n\text{ and }j+m>m\right)  .
\end{align*}
This rewrites as $d_{i,j}=d_{i,j}^{\prime}$ (since $\left(  i+n\right)  -n=i$
and $\left(  j+m\right)  -m=j$). Qed.}. In other words, $\left(
d_{i,j}\right)  _{1\leq i\leq n^{\prime},\ 1\leq j\leq m^{\prime}}=\left(
d_{i,j}^{\prime}\right)  _{1\leq i\leq n^{\prime},\ 1\leq j\leq m^{\prime}}$.
Now,%
\[
D=\left(  d_{i,j}\right)  _{1\leq i\leq n^{\prime},\ 1\leq j\leq m^{\prime}%
}=\left(  d_{i,j}^{\prime}\right)  _{1\leq i\leq n^{\prime},\ 1\leq j\leq
m^{\prime}}=D^{\prime}.
\]

We now have shown that $A=A^{\prime}$, $B=B^{\prime}$, $C=C^{\prime}$ and
$D=D^{\prime}$. This proves Lemma \ref{lem.sol.block2x2.jacobi.equal}.
\end{proof}
\end{verlong}

\begin{proof}
[Solution to Exercise \ref{exe.block2x2.jacobi}.]The matrix $\left(
\begin{array}
[c]{cc}%
A^{\prime} & B^{\prime}\\
C^{\prime} & D^{\prime}%
\end{array}
\right)  $ is the inverse of the matrix $\left(
\begin{array}
[c]{cc}%
A & B\\
C & D
\end{array}
\right)  $. In other words, we have%
\begin{align*}
\left(
\begin{array}
[c]{cc}%
A & B\\
C & D
\end{array}
\right)  \left(
\begin{array}
[c]{cc}%
A^{\prime} & B^{\prime}\\
C^{\prime} & D^{\prime}%
\end{array}
\right)   &  =I_{n+m}\ \ \ \ \ \ \ \ \ \ \text{and}\\
\left(
\begin{array}
[c]{cc}%
A^{\prime} & B^{\prime}\\
C^{\prime} & D^{\prime}%
\end{array}
\right)  \left(
\begin{array}
[c]{cc}%
A & B\\
C & D
\end{array}
\right)   &  =I_{n+m}.
\end{align*}
But Lemma \ref{lem.sol.block2x2.jacobi.InIm} yields
\[
\left(
\begin{array}
[c]{cc}%
I_{n} & 0_{n\times m}\\
0_{m\times n} & I_{m}%
\end{array}
\right)  =I_{n+m}=\left(
\begin{array}
[c]{cc}%
A^{\prime} & B^{\prime}\\
C^{\prime} & D^{\prime}%
\end{array}
\right)  \left(
\begin{array}
[c]{cc}%
A & B\\
C & D
\end{array}
\right)  =\left(
\begin{array}
[c]{cc}%
A^{\prime}A+B^{\prime}C & A^{\prime}B+B^{\prime}D\\
C^{\prime}A+D^{\prime}C & C^{\prime}B+D^{\prime}D
\end{array}
\right)
\]
(by Exercise \ref{exe.block2x2.mult} (applied to $n$, $m$, $n$, $m$, $n$, $m$,
$A^{\prime}$, $B^{\prime}$, $C^{\prime}$, $D^{\prime}$, $A$, $B$, $C$ and $D$
instead of $n$, $n^{\prime}$, $m$, $m^{\prime}$, $\ell$, $\ell^{\prime}$, $A$,
$B$, $C$, $D$, $A^{\prime}$, $B^{\prime}$, $C^{\prime}$ and $D^{\prime}$)).
Thus, Lemma \ref{lem.sol.block2x2.jacobi.equal} (applied to $n$, $m$, $n$,
$m$, $I_{n}$, $A^{\prime}A+B^{\prime}C$, $0_{n\times m}$, $A^{\prime
}B+B^{\prime}D$, $0_{m\times n}$, $C^{\prime}A+D^{\prime}C$, $I_{m}$ and
$C^{\prime}B+D^{\prime}D$ instead of $n$, $n^{\prime}$, $m$, $m^{\prime}$,
$A$, $A^{\prime}$, $B$, $B^{\prime}$, $C$, $C^{\prime}$, $D$ and $D^{\prime}$)
yields that $I_{n}=A^{\prime}A+B^{\prime}C$, $0_{n\times m}=A^{\prime
}B+B^{\prime}D$, $0_{m\times n}=C^{\prime}A+D^{\prime}C$ and $I_{m}=C^{\prime
}B+D^{\prime}D$.

Recall that $\det\left(  I_{n}\right)  =1$. The same argument (applied to $m$
instead of $n$) yields $\det\left(  I_{m}\right)  =1$. But from $I_{m}%
=C^{\prime}B+D^{\prime}D$, we obtain $C^{\prime}B+D^{\prime}D=I_{m}$ and thus%
\[
\det\underbrace{\left(  C^{\prime}B+D^{\prime}D\right)  }_{=I_{m}}=\det\left(
I_{m}\right)  =1.
\]

From $0_{m\times n}=C^{\prime}A+D^{\prime}C$, we obtain $C^{\prime
}A=-D^{\prime}C$. Hence, Exercise \ref{exe.block2x2.VB+WD} (applied to
$W=D^{\prime}$ and $V=C^{\prime}$) yields
\[
\det\left(  D^{\prime}\right)  \cdot\det\left(
\begin{array}
[c]{cc}%
A & B\\
C & D
\end{array}
\right)  =\det A\cdot\underbrace{\det\left(  C^{\prime}B+D^{\prime}D\right)
}_{=1}=\det A.
\]
Thus,%
\[
\det A=\det\left(  D^{\prime}\right)  \cdot\det\left(
\begin{array}
[c]{cc}%
A & B\\
C & D
\end{array}
\right)  =\det\left(
\begin{array}
[c]{cc}%
A & B\\
C & D
\end{array}
\right)  \cdot\det\left(  D^{\prime}\right)  .
\]
This solves Exercise \ref{exe.block2x2.jacobi}.
\end{proof}

\subsection{Solution to Exercise \ref{exe.block2x2.jacobi.rewr}}

Our solution to Exercise \ref{exe.block2x2.jacobi.rewr} below will use the
following simple lemma:

\begin{lemma}
\label{lem.sol.block2x2.jacobi.rewr.1}Let $n\in\mathbb{N}$ and $m\in
\mathbb{N}$. Let $A\in\mathbb{K}^{n\times m}$. Let $k\in\left\{
0,1,\ldots,n\right\}  $ and $\ell\in\left\{  0,1,\ldots,m\right\}  $. Then,%
\[
A=\left(
\begin{array}
[c]{cc}%
\operatorname*{sub}\nolimits_{1,2,\ldots,k}^{1,2,\ldots,\ell}A &
\operatorname*{sub}\nolimits_{1,2,\ldots,k}^{\ell+1,\ell+2,\ldots,m}A\\
\operatorname*{sub}\nolimits_{k+1,k+2,\ldots,n}^{1,2,\ldots,\ell}A &
\operatorname*{sub}\nolimits_{k+1,k+2,\ldots,n}^{\ell+1,\ell+2,\ldots,m}A
\end{array}
\right)  .
\]

\end{lemma}

\begin{vershort}
\begin{proof}
[Proof of Lemma \ref{lem.sol.block2x2.jacobi.rewr.1}.]In visual language,
Lemma \ref{lem.sol.block2x2.jacobi.rewr.1} states a triviality: It says that
if we cut the matrix $A$ horizontally (between its $k$-th and $\left(
k+1\right)  $-st rows) and vertically (between its $\ell$-th and $\left(
\ell+1\right)  $-st columns), then we obtain four little matrices (namely,
$\operatorname*{sub}\nolimits_{1,2,\ldots,k}^{1,2,\ldots,\ell}A$,
$\operatorname*{sub}\nolimits_{1,2,\ldots,k}^{\ell+1,\ell+2,\ldots,m}A$,
$\operatorname*{sub}\nolimits_{k+1,k+2,\ldots,n}^{1,2,\ldots,\ell}A$ and
$\operatorname*{sub}\nolimits_{k+1,k+2,\ldots,n}^{\ell+1,\ell+2,\ldots,m}A$)
which can be assembled back to form $A$ (using the block-matrix construction).
Turning this into a formal proof is straightforward.
\end{proof}
\end{vershort}

\begin{verlong}
\begin{proof}
[Proof of Lemma \ref{lem.sol.block2x2.jacobi.rewr.1}.]Write the matrix
$A\in\mathbb{K}^{n\times m}$ in the form $A=\left(  a_{i,j}\right)  _{1\leq
i\leq n,\ 1\leq j\leq m}$.

We have $k\in\left\{  0,1,\ldots,n\right\}  $, thus $k\in\mathbb{N}$ and
$k\leq n$. From $k\leq n$, we obtain $n-k\geq0$ and thus $n-k\in\mathbb{N}$.

We have $\ell\in\left\{  0,1,\ldots,m\right\}  $, thus $\ell\in\mathbb{N}$ and
$\ell\leq m$. From $\ell\leq m$, we obtain $m-\ell\geq0$ and thus $m-\ell
\in\mathbb{N}$.

The definition of $\operatorname*{sub}\nolimits_{1,2,\ldots,k}^{1,2,\ldots
,\ell}A$ yields $\operatorname*{sub}\nolimits_{1,2,\ldots,k}^{1,2,\ldots,\ell
}A=\left(  a_{i,j}\right)  _{1\leq i\leq k,\ 1\leq j\leq\ell}$ (since
$A=\left(  a_{i,j}\right)  _{1\leq i\leq n,\ 1\leq j\leq m}$).

The definition of $\operatorname*{sub}\nolimits_{1,2,\ldots,k}^{\ell
+1,\ell+2,\ldots,m}A$ yields $\operatorname*{sub}\nolimits_{1,2,\ldots
,k}^{\ell+1,\ell+2,\ldots,m}A=\left(  a_{i,\ell+j}\right)  _{1\leq i\leq
k,\ 1\leq j\leq m-\ell}$ (since $A=\left(  a_{i,j}\right)  _{1\leq i\leq
n,\ 1\leq j\leq m}$).

The definition of $\operatorname*{sub}\nolimits_{k+1,k+2,\ldots,n}%
^{1,2,\ldots,\ell}A$ yields $\operatorname*{sub}\nolimits_{k+1,k+2,\ldots
,n}^{1,2,\ldots,\ell}A=\left(  a_{k+i,j}\right)  _{1\leq i\leq n-k,\ 1\leq
j\leq\ell}$ (since $A=\left(  a_{i,j}\right)  _{1\leq i\leq n,\ 1\leq j\leq
m}$).

The definition of $\operatorname*{sub}\nolimits_{k+1,k+2,\ldots,n}%
^{\ell+1,\ell+2,\ldots,m}A$ yields $\operatorname*{sub}%
\nolimits_{k+1,k+2,\ldots,n}^{\ell+1,\ell+2,\ldots,m}A=\left(  a_{k+i,\ell
+j}\right)  _{1\leq i\leq n-k,\ 1\leq j\leq m-\ell}$ (since $A=\left(
a_{i,j}\right)  _{1\leq i\leq n,\ 1\leq j\leq m}$).

We have $\operatorname*{sub}\nolimits_{1,2,\ldots,k}^{1,2,\ldots,\ell
}A=\left(  a_{i,j}\right)  _{1\leq i\leq k,\ 1\leq j\leq\ell}$,
$\operatorname*{sub}\nolimits_{1,2,\ldots,k}^{\ell+1,\ell+2,\ldots,m}A=\left(
a_{i,\ell+j}\right)  _{1\leq i\leq k,\ 1\leq j\leq m-\ell}$,
$\operatorname*{sub}\nolimits_{k+1,k+2,\ldots,n}^{1,2,\ldots,\ell}A=\left(
a_{k+i,j}\right)  _{1\leq i\leq n-k,\ 1\leq j\leq\ell}$ and
$\operatorname*{sub}\nolimits_{k+1,k+2,\ldots,n}^{\ell+1,\ell+2,\ldots
,m}A=\left(  a_{k+i,\ell+j}\right)  _{1\leq i\leq n-k,\ 1\leq j\leq m-\ell}$.
Hence, the definition of $\left(
\begin{array}
[c]{cc}%
\operatorname*{sub}\nolimits_{1,2,\ldots,k}^{1,2,\ldots,\ell}A &
\operatorname*{sub}\nolimits_{1,2,\ldots,k}^{\ell+1,\ell+2,\ldots,m}A\\
\operatorname*{sub}\nolimits_{k+1,k+2,\ldots,n}^{1,2,\ldots,\ell}A &
\operatorname*{sub}\nolimits_{k+1,k+2,\ldots,n}^{\ell+1,\ell+2,\ldots,m}A
\end{array}
\right)  $ yields%
\begin{align}
&  \left(
\begin{array}
[c]{cc}%
\operatorname*{sub}\nolimits_{1,2,\ldots,k}^{1,2,\ldots,\ell}A &
\operatorname*{sub}\nolimits_{1,2,\ldots,k}^{\ell+1,\ell+2,\ldots,m}A\\
\operatorname*{sub}\nolimits_{k+1,k+2,\ldots,n}^{1,2,\ldots,\ell}A &
\operatorname*{sub}\nolimits_{k+1,k+2,\ldots,n}^{\ell+1,\ell+2,\ldots,m}A
\end{array}
\right) \nonumber\\
&  =\left(
\begin{cases}
a_{i,j}, & \text{if }i\leq k\text{ and }j\leq\ell;\\
a_{i,\ell+\left(  j-\ell\right)  }, & \text{if }i\leq k\text{ and }j>\ell;\\
a_{k+\left(  i-k\right)  ,j}, & \text{if }i>k\text{ and }j\leq\ell;\\
a_{k+\left(  i-k\right)  ,\ell+\left(  j-\ell\right)  }, & \text{if }i>k\text{
and }j>\ell
\end{cases}
\right)  _{1\leq i\leq k+\left(  n-k\right)  ,\ 1\leq j\leq\ell+\left(
m-\ell\right)  }\nonumber\\
&  =\left(
\begin{cases}
a_{i,j}, & \text{if }i\leq k\text{ and }j\leq\ell;\\
a_{i,\ell+\left(  j-\ell\right)  }, & \text{if }i\leq k\text{ and }j>\ell;\\
a_{k+\left(  i-k\right)  ,j}, & \text{if }i>k\text{ and }j\leq\ell;\\
a_{k+\left(  i-k\right)  ,\ell+\left(  j-\ell\right)  }, & \text{if }i>k\text{
and }j>\ell
\end{cases}
\right)  _{1\leq i\leq n,\ 1\leq j\leq m}
\label{pf.lem.sol.block2x2.jacobi.rewr.1.1}%
\end{align}
(since $k+\left(  n-k\right)  =n$ and $\ell+\left(  m-\ell\right)  =m$).

But every $\left(  i,j\right)  \in\left\{  1,2,\ldots,n\right\}
\times\left\{  1,2,\ldots,m\right\}  $ satisfies%
\begin{align}
&
\begin{cases}
a_{i,j}, & \text{if }i\leq k\text{ and }j\leq\ell;\\
a_{i,\ell+\left(  j-\ell\right)  }, & \text{if }i\leq k\text{ and }j>\ell;\\
a_{k+\left(  i-k\right)  ,j}, & \text{if }i>k\text{ and }j\leq\ell;\\
a_{k+\left(  i-k\right)  ,\ell+\left(  j-\ell\right)  }, & \text{if }i>k\text{
and }j>\ell
\end{cases}
\nonumber\\
&  =%
\begin{cases}
a_{i,j}, & \text{if }i\leq k\text{ and }j\leq\ell;\\
a_{i,j}, & \text{if }i\leq k\text{ and }j>\ell;\\
a_{k+\left(  i-k\right)  ,j}, & \text{if }i>k\text{ and }j\leq\ell;\\
a_{k+\left(  i-k\right)  ,j}, & \text{if }i>k\text{ and }j>\ell
\end{cases}
\ \ \ \ \ \ \ \ \ \ \left(  \text{since }\ell+\left(  j-\ell\right)  =j\right)
\nonumber\\
&  =%
\begin{cases}
a_{i,j}, & \text{if }i\leq k\text{ and }j\leq\ell;\\
a_{i,j}, & \text{if }i\leq k\text{ and }j>\ell;\\
a_{i,j}, & \text{if }i>k\text{ and }j\leq\ell;\\
a_{i,j}, & \text{if }i>k\text{ and }j>\ell
\end{cases}
\ \ \ \ \ \ \ \ \ \ \left(  \text{since }k+\left(  i-k\right)  =i\right)
\nonumber\\
&  =a_{i,j}. \label{pf.lem.sol.block2x2.jacobi.rewr.1.2}%
\end{align}
Now, (\ref{pf.lem.sol.block2x2.jacobi.rewr.1.1}) becomes%
\begin{align*}
&  \left(
\begin{array}
[c]{cc}%
\operatorname*{sub}\nolimits_{1,2,\ldots,k}^{1,2,\ldots,\ell}A &
\operatorname*{sub}\nolimits_{1,2,\ldots,k}^{\ell+1,\ell+2,\ldots,m}A\\
\operatorname*{sub}\nolimits_{k+1,k+2,\ldots,n}^{1,2,\ldots,\ell}A &
\operatorname*{sub}\nolimits_{k+1,k+2,\ldots,n}^{\ell+1,\ell+2,\ldots,m}A
\end{array}
\right) \\
&  =\left(  \underbrace{%
\begin{cases}
a_{i,j}, & \text{if }i\leq k\text{ and }j\leq\ell;\\
a_{i,\ell+\left(  j-\ell\right)  }, & \text{if }i\leq k\text{ and }j>\ell;\\
a_{k+\left(  i-k\right)  ,j}, & \text{if }i>k\text{ and }j\leq\ell;\\
a_{k+\left(  i-k\right)  ,\ell+\left(  j-\ell\right)  }, & \text{if }i>k\text{
and }j>\ell
\end{cases}
}_{\substack{=a_{i,j}\\\text{(by (\ref{pf.lem.sol.block2x2.jacobi.rewr.1.2}%
))}}}\right)  _{1\leq i\leq n,\ 1\leq j\leq m}\\
&  =\left(  a_{i,j}\right)  _{1\leq i\leq n,\ 1\leq j\leq m}=A.
\end{align*}
This proves Lemma \ref{lem.sol.block2x2.jacobi.rewr.1}.
\end{proof}
\end{verlong}

\begin{proof}
[Solution to Exercise \ref{exe.block2x2.jacobi.rewr}.]Lemma
\ref{lem.sol.block2x2.jacobi.rewr.1} (applied to $m=n$ and $\ell=k$) yields%
\begin{equation}
A=\left(
\begin{array}
[c]{cc}%
\operatorname*{sub}\nolimits_{1,2,\ldots,k}^{1,2,\ldots,k}A &
\operatorname*{sub}\nolimits_{1,2,\ldots,k}^{k+1,k+2,\ldots,n}A\\
\operatorname*{sub}\nolimits_{k+1,k+2,\ldots,n}^{1,2,\ldots,k}A &
\operatorname*{sub}\nolimits_{k+1,k+2,\ldots,n}^{k+1,k+2,\ldots,n}A
\end{array}
\right)  . \label{sol.block2x2.jacobi.rewr.A=}%
\end{equation}
But the matrix $A\in\mathbb{K}^{n\times n}$ is invertible. Its inverse is
$A^{-1}\in\mathbb{K}^{n\times n}$. Lemma \ref{lem.sol.block2x2.jacobi.rewr.1}
(applied to $n$, $k$ and $A^{-1}$ instead of $m$, $\ell$ and $A$) yields%
\begin{equation}
A^{-1}=\left(
\begin{array}
[c]{cc}%
\operatorname*{sub}\nolimits_{1,2,\ldots,k}^{1,2,\ldots,k}\left(
A^{-1}\right)  & \operatorname*{sub}\nolimits_{1,2,\ldots,k}^{k+1,k+2,\ldots
,n}\left(  A^{-1}\right) \\
\operatorname*{sub}\nolimits_{k+1,k+2,\ldots,n}^{1,2,\ldots,k}\left(
A^{-1}\right)  & \operatorname*{sub}\nolimits_{k+1,k+2,\ldots,n}%
^{k+1,k+2,\ldots,n}\left(  A^{-1}\right)
\end{array}
\right)  . \label{sol.block2x2.jacobi.rewr.A-1=}%
\end{equation}

\begin{verlong}
We have $k\in\left\{  0,1,\ldots,n\right\}  $, thus $k\in\mathbb{N}$ and
$k\leq n$. From $k\leq n$, we obtain $n-k\geq0$ and thus $n-k\in\mathbb{N}$.
\end{verlong}

Now, recall that the matrix $A$ is invertible, and its inverse is $A^{-1}$. In
view of the equalities (\ref{sol.block2x2.jacobi.rewr.A=}) and
(\ref{sol.block2x2.jacobi.rewr.A-1=}), this rewrites as follows: The matrix
\newline$\left(
\begin{array}
[c]{cc}%
\operatorname*{sub}\nolimits_{1,2,\ldots,k}^{1,2,\ldots,k}A &
\operatorname*{sub}\nolimits_{1,2,\ldots,k}^{k+1,k+2,\ldots,n}A\\
\operatorname*{sub}\nolimits_{k+1,k+2,\ldots,n}^{1,2,\ldots,k}A &
\operatorname*{sub}\nolimits_{k+1,k+2,\ldots,n}^{k+1,k+2,\ldots,n}A
\end{array}
\right)  $ is invertible, and its inverse is \newline$\left(
\begin{array}
[c]{cc}%
\operatorname*{sub}\nolimits_{1,2,\ldots,k}^{1,2,\ldots,k}\left(
A^{-1}\right)  & \operatorname*{sub}\nolimits_{1,2,\ldots,k}^{k+1,k+2,\ldots
,n}\left(  A^{-1}\right) \\
\operatorname*{sub}\nolimits_{k+1,k+2,\ldots,n}^{1,2,\ldots,k}\left(
A^{-1}\right)  & \operatorname*{sub}\nolimits_{k+1,k+2,\ldots,n}%
^{k+1,k+2,\ldots,n}\left(  A^{-1}\right)
\end{array}
\right)  $. Hence, Exercise \ref{exe.block2x2.jacobi} (applied to $k$, $n-k$,
$\operatorname*{sub}\nolimits_{1,2,\ldots,k}^{1,2,\ldots,k}A$,
$\operatorname*{sub}\nolimits_{1,2,\ldots,k}^{k+1,k+2,\ldots,n}A$,
$\operatorname*{sub}\nolimits_{k+1,k+2,\ldots,n}^{1,2,\ldots,k}A$,
$\operatorname*{sub}\nolimits_{k+1,k+2,\ldots,n}^{k+1,k+2,\ldots,n}A$,
$\operatorname*{sub}\nolimits_{1,2,\ldots,k}^{1,2,\ldots,k}\left(
A^{-1}\right)  $, $\operatorname*{sub}\nolimits_{1,2,\ldots,k}^{k+1,k+2,\ldots
,n}\left(  A^{-1}\right)  $, $\operatorname*{sub}\nolimits_{k+1,k+2,\ldots
,n}^{1,2,\ldots,k}\left(  A^{-1}\right)  $ and $\operatorname*{sub}%
\nolimits_{k+1,k+2,\ldots,n}^{k+1,k+2,\ldots,n}\left(  A^{-1}\right)  $
instead of $n$, $m$, $A$, $B$, $C$, $D$, $A^{\prime}$, $B^{\prime}$,
$C^{\prime}$ and $D^{\prime}$) yields%
\begin{align*}
&  \det\left(  \operatorname*{sub}\nolimits_{1,2,\ldots,k}^{1,2,\ldots
,k}A\right) \\
&  =\det\underbrace{\left(
\begin{array}
[c]{cc}%
\operatorname*{sub}\nolimits_{1,2,\ldots,k}^{1,2,\ldots,k}A &
\operatorname*{sub}\nolimits_{1,2,\ldots,k}^{k+1,k+2,\ldots,n}A\\
\operatorname*{sub}\nolimits_{k+1,k+2,\ldots,n}^{1,2,\ldots,k}A &
\operatorname*{sub}\nolimits_{k+1,k+2,\ldots,n}^{k+1,k+2,\ldots,n}A
\end{array}
\right)  }_{\substack{=A\\\text{(by (\ref{sol.block2x2.jacobi.rewr.A=}))}%
}}\cdot\det\left(  \operatorname*{sub}\nolimits_{k+1,k+2,\ldots,n}%
^{k+1,k+2,\ldots,n}\left(  A^{-1}\right)  \right) \\
&  =\det A\cdot\det\left(  \operatorname*{sub}\nolimits_{k+1,k+2,\ldots
,n}^{k+1,k+2,\ldots,n}\left(  A^{-1}\right)  \right)  .
\end{align*}
This solves Exercise \ref{exe.block2x2.jacobi.rewr}.
\end{proof}

\subsection{Solution to Exercise \ref{exe.unrows.basics}}

\begin{vershort}
\begin{proof}
[Solution to Exercise \ref{exe.unrows.basics}.]\textbf{(a)} All claims of
Proposition \ref{prop.unrows.basics} and Proposition
\ref{prop.unrows.basics-I} are \textquotedblleft clear by
inspection\textquotedblright, in the sense that the reader will have no
difficulty convincing themselves of their validity just by studying an example
and seeing \textquotedblleft what is going on\textquotedblright. Let me give
some more detailed proofs. (The following proofs are not very formal; but I
will outline a formal proof of Proposition \ref{prop.unrows.basics}
\textbf{(h)} in a footnote, and I trust that the reader can use the same
methods to formally prove the rest of Proposition \ref{prop.unrows.basics} if
so inclined.)

\begin{proof}
[Proof of Proposition \ref{prop.unrows.basics}.]Proposition
\ref{prop.unrows.basics} \textbf{(a)} follows from the definitions of
$A_{u,\bullet}$ and of $\operatorname*{rows}\nolimits_{u}A$ (indeed, these
definitions show that both $A_{u,\bullet}$ and $\operatorname*{rows}%
\nolimits_{u}A$ are the $u$-th row of $A$). Similarly, Proposition
\ref{prop.unrows.basics} \textbf{(b)} follows from the definitions of
$A_{\bullet,v}$ and $\operatorname*{cols}\nolimits_{v}A$.

Proposition \ref{prop.unrows.basics} \textbf{(c)} is obvious\footnote{In fact:
\par
\begin{itemize}
\item The matrix $\left(  A_{\bullet,\sim v}\right)  _{\sim u,\bullet}$ is the
matrix obtained from $A$ by first removing the $v$-th column and then removing
the $u$-th row.
\par
\item The matrix $\left(  A_{\sim u,\bullet}\right)  _{\bullet,\sim v}$ is the
matrix obtained from $A$ by first removing the $u$-th row and then removing
the $v$-th column.
\par
\item The matrix $A_{\sim u,\sim v}$ is the matrix obtained from $A$ by
removing the $u$-th row and the $v$-th column at the same time.
\end{itemize}
\par
Thus it is clear that all three of these matrices are equal.}.

Proposition \ref{prop.unrows.basics} \textbf{(d)} claims that the $w$-th
column of the matrix $A_{\bullet,\sim v}$ equals the $w$-th column of the
matrix $A$ (whenever $v\in\left\{  1,2,\ldots,m\right\}  $ and $w\in\left\{
1,2,\ldots,v-1\right\}  $). Let us prove this: Let $v\in\left\{
1,2,\ldots,m\right\}  $ and $w\in\left\{  1,2,\ldots,v-1\right\}  $. Thus,
$w<v$ (since $w\in\left\{  1,2,\ldots,v-1\right\}  $). But the matrix
$A_{\bullet,\sim v}$ results from the matrix $A$ by removing the $v$-th
column; clearly, this removal does not change the $w$-th column (because
$w<v$). Thus, the $w$-th column of the matrix $A_{\bullet,\sim v}$ equals the
$w$-th column of the matrix $A$. So Proposition \ref{prop.unrows.basics}
\textbf{(d)} is proven.\footnote{If you found this proof insufficiently
rigorous, let me show a formal proof of Proposition \ref{prop.unrows.basics}
\textbf{(d)}. First, I shall introduce a notation:
\par
\begin{itemize}
\item For every $j\in\mathbb{Z}$ and $r\in\mathbb{Z}$, let $\mathbf{d}%
_{r}\left(  j\right)  $ be the integer $%
\begin{cases}
j, & \text{if }j<r;\\
j+1, & \text{if }j\geq r
\end{cases}
$.
\end{itemize}
\par
We make the following observation:
\par
\textit{Observation 1:} If $n\in\mathbb{N}$ and $m\in\mathbb{N}$, if
$A=\left(  a_{i,j}\right)  _{1\leq i\leq n,\ 1\leq j\leq m}\in\mathbb{K}%
^{n\times m}$ is an $n\times m$-matrix, and if $v$ is an element of $\left\{
1,2,\ldots,m\right\}  $, then $A_{\bullet,\sim v}=\left(  a_{i,\mathbf{d}%
_{v}\left(  j\right)  }\right)  _{1\leq i\leq n,\ 1\leq j\leq m-1}$.
\par
\textit{Proof of Observation 1:} Let $n\in\mathbb{N}$ and $m\in\mathbb{N}$.
Let $A=\left(  a_{i,j}\right)  _{1\leq i\leq n,\ 1\leq j\leq m}\in
\mathbb{K}^{n\times m}$ be an $n\times m$-matrix. Let $v$ be an element of
$\left\{  1,2,\ldots,m\right\}  $. Recalling how $\mathbf{d}_{v}\left(
j\right)  $ is defined for every $j\in\mathbb{Z}$, we see that%
\begin{align*}
\left(  \mathbf{d}_{v}\left(  1\right)  ,\mathbf{d}_{v}\left(  2\right)
,\ldots,\mathbf{d}_{v}\left(  m-1\right)  \right)   &  =\left(  1,2,\ldots
,v-1,v+1,v+2,\ldots,m\right) \\
&  =\left(  1,2,\ldots,\widehat{v},\ldots,m\right)  .
\end{align*}
But the definition of $A_{\bullet,\sim v}$ yields
\begin{align*}
A_{\bullet,\sim v}  &  =\operatorname*{cols}\nolimits_{1,2,\ldots
,\widehat{v},\ldots,m}A=\operatorname*{cols}\nolimits_{\mathbf{d}_{v}\left(
1\right)  ,\mathbf{d}_{v}\left(  2\right)  ,\ldots,\mathbf{d}_{v}\left(
m-1\right)  }A\\
&  \ \ \ \ \ \ \ \ \ \ \left(  \text{since }\left(  1,2,\ldots,\widehat{v}%
,\ldots,m\right)  =\left(  \mathbf{d}_{v}\left(  1\right)  ,\mathbf{d}%
_{v}\left(  2\right)  ,\ldots,\mathbf{d}_{v}\left(  m-1\right)  \right)
\right) \\
&  =\left(  a_{i,\mathbf{d}_{v}\left(  y\right)  }\right)  _{1\leq i\leq
n,\ 1\leq y\leq m-1}\ \ \ \ \ \ \ \ \ \ \left(
\begin{array}
[c]{c}%
\text{by the definition of }\operatorname*{cols}\nolimits_{\mathbf{d}%
_{v}\left(  1\right)  ,\mathbf{d}_{v}\left(  2\right)  ,\ldots,\mathbf{d}%
_{v}\left(  m-1\right)  }A\\
\text{(since }A=\left(  a_{i,j}\right)  _{1\leq i\leq n,\ 1\leq j\leq
m}\text{)}%
\end{array}
\right) \\
&  =\left(  a_{i,\mathbf{d}_{v}\left(  j\right)  }\right)  _{1\leq i\leq
n,\ 1\leq j\leq m-1}\ \ \ \ \ \ \ \ \ \ \left(  \text{here, we have renamed
the index }\left(  i,y\right)  \text{ as }\left(  i,j\right)  \right)  .
\end{align*}
This proves Observation 1.
\par
Let us now prove Proposition \ref{prop.unrows.basics} \textbf{(d)}. Indeed,
let $v\in\left\{  1,2,\ldots,m\right\}  $ and $w\in\left\{  1,2,\ldots
,v-1\right\}  $. Then, $w<v$ (since $w\in\left\{  1,2,\ldots,v-1\right\}  $).
Write the $n\times m$-matrix $A$ as $A=\left(  a_{i,j}\right)  _{1\leq i\leq
n,\ 1\leq j\leq m}$. The definition of $\mathbf{d}_{v}\left(  w\right)  $ now
yields $\mathbf{d}_{v}\left(  w\right)  =%
\begin{cases}
w, & \text{if }w<v;\\
w+1, & \text{if }w\geq v
\end{cases}
=w$ (since $w<v$). Now, the definition of $\left(  A_{\bullet,\sim v}\right)
_{\bullet,w}$ yields%
\begin{align*}
\left(  A_{\bullet,\sim v}\right)  _{\bullet,w}  &  =\left(  \text{the
}w\text{-th column of the matrix }A_{\bullet,\sim v}\right) \\
&  =\left(
\begin{array}
[c]{c}%
a_{1,\mathbf{d}_{v}\left(  w\right)  }\\
a_{2,\mathbf{d}_{v}\left(  w\right)  }\\
\vdots\\
a_{n,\mathbf{d}_{v}\left(  w\right)  }%
\end{array}
\right)  \ \ \ \ \ \ \ \ \ \ \left(  \text{since Observation 1 yields
}A_{\bullet,\sim v}=\left(  a_{i,\mathbf{d}_{v}\left(  j\right)  }\right)
_{1\leq i\leq n,\ 1\leq j\leq m-1}\right) \\
&  =\left(
\begin{array}
[c]{c}%
a_{1,w}\\
a_{2,w}\\
\vdots\\
a_{n,w}%
\end{array}
\right)  \ \ \ \ \ \ \ \ \ \ \left(  \text{since }\mathbf{d}_{v}\left(
w\right)  =w\right) \\
&  =\left(  \text{the }w\text{-th column of the matrix }A\right)
\ \ \ \ \ \ \ \ \ \ \left(  \text{since }A=\left(  a_{i,j}\right)  _{1\leq
i\leq n,\ 1\leq j\leq m}\right) \\
&  =A_{\bullet,w}.
\end{align*}
This proves Proposition \ref{prop.unrows.basics} \textbf{(d)}.}

Proposition \ref{prop.unrows.basics} \textbf{(e)} claims that the $w$-th
column of the matrix $A_{\bullet,\sim v}$ equals the $\left(  w+1\right)  $-th
column of the matrix $A$ (whenever $v\in\left\{  1,2,\ldots,m\right\}  $ and
$w\in\left\{  v,v+1,\ldots,m-1\right\}  $). Let us prove this:

Let $v\in\left\{  1,2,\ldots,m\right\}  $ and $w\in\left\{  v,v+1,\ldots
,m-1\right\}  $. Thus, $w+1>w\geq v$ (since $w\in\left\{  v,v+1,\ldots
,m-1\right\}  $). But the matrix $A_{\bullet,\sim v}$ results from the matrix
$A$ by removing the $v$-th column; this removal moves the $\left(  w+1\right)
$-th column of $A$ one step leftwards (since $w+1>v$). In other words, the
$\left(  w+1\right)  $-th column of $A$ becomes the $w$-th column of the new
matrix $A_{\bullet,\sim v}$. In other words, the $w$-th column of the matrix
$A_{\bullet,\sim v}$ equals the $\left(  w+1\right)  $-th column of the matrix
$A$. So Proposition \ref{prop.unrows.basics} \textbf{(e)} is
proven.\footnote{Again, a more rigorous proof of Proposition
\ref{prop.unrows.basics} \textbf{(e)} can be obtained similarly to our
rigorous proof of Proposition \ref{prop.unrows.basics} \textbf{(d)} in the
previous footnote. (The main difference is that we have $\mathbf{d}_{v}\left(
w\right)  =w+1$ instead of $\mathbf{d}_{v}\left(  w\right)  =w$ this time.)}

Proposition \ref{prop.unrows.basics} \textbf{(f)} is the analogue of
Proposition \ref{prop.unrows.basics} \textbf{(d)} for rows instead of columns
(with $v$ renamed as $u$); its proof is equally analogous.

Proposition \ref{prop.unrows.basics} \textbf{(g)} is the analogue of
Proposition \ref{prop.unrows.basics} \textbf{(e)} for rows instead of columns
(with $v$ renamed as $u$); its proof is equally analogous.

Let us now prove Proposition \ref{prop.unrows.basics} \textbf{(h)}: Let
$v\in\left\{  1,2,\ldots,m\right\}  $ and $w\in\left\{  1,2,\ldots
,v-1\right\}  $. We need to show that $\left(  A_{\bullet,\sim v}\right)
_{\bullet,\sim w}=\operatorname*{cols}\nolimits_{1,2,\ldots,\widehat{w}%
,\ldots,\widehat{v},\ldots,m}A$. The matrix $\left(  A_{\bullet,\sim
v}\right)  _{\bullet,\sim w}$ is obtained from the matrix $A$ by first
removing the $v$-th column (thus obtaining an $n\times\left(  m-1\right)
$-matrix) and then removing the $w$-th column (from the resulting
$n\times\left(  m-1\right)  $-matrix). Since $w<v$ (because $w\in\left\{
1,2,\ldots,v-1\right\}  $), the removal of the $v$-th column of $A$ did not
affect the $w$-th column, and therefore the two successive removals could be
replaced by a simultaneous removal of both the $w$-th and the $v$-th columns
from the matrix $A$. In other words, the matrix obtained from the matrix $A$
by first removing the $v$-th column and then removing the $w$-th column is
identical with the matrix obtained from $A$ by simultaneously removing both
the $w$-th and the $v$-th columns. But the former matrix is $\left(
A_{\bullet,\sim v}\right)  _{\bullet,\sim w}$, whereas the latter matrix is
$\operatorname*{cols}\nolimits_{1,2,\ldots,\widehat{w},\ldots,\widehat{v}%
,\ldots,m}A$. Hence, the identity of these two matrices rewrites as $\left(
A_{\bullet,\sim v}\right)  _{\bullet,\sim w}=\operatorname*{cols}%
\nolimits_{1,2,\ldots,\widehat{w},\ldots,\widehat{v},\ldots,m}A$. This proves
Proposition \ref{prop.unrows.basics} \textbf{(h)}.\footnote{If you found this
proof insufficiently rigorous, let me show a formal proof of Proposition
\ref{prop.unrows.basics} \textbf{(h)}. First, we have $r\leq s-1$ (since $r<s$
and since $r$ and $s$ are integers). Now, I shall introduce two notations:
\par
\begin{itemize}
\item For every $j\in\mathbb{Z}$ and $r\in\mathbb{Z}$, let $\mathbf{d}%
_{r}\left(  j\right)  $ be the integer $%
\begin{cases}
j, & \text{if }j<r;\\
j+1, & \text{if }j\geq r
\end{cases}
$.
\par
\item For every $j\in\mathbb{Z}$, $r\in\mathbb{Z}$ and $s\in\mathbb{Z}$
satisfying $r<s$, we let $\mathbf{d}_{r,s}\left(  j\right)  $ be the integer $%
\begin{cases}
j, & \text{if }j<r;\\
j+1, & \text{if }r\leq j<s-1;\\
j+2, & \text{if }s-1\leq j
\end{cases}
$. (This is well-defined, because $r\leq s-1$.)
\end{itemize}
\par
We make the following three observations:
\par
\textit{Observation 1:} If $n\in\mathbb{N}$ and $m\in\mathbb{N}$, if
$A=\left(  a_{i,j}\right)  _{1\leq i\leq n,\ 1\leq j\leq m}\in\mathbb{K}%
^{n\times m}$ is an $n\times m$-matrix, and if $v$ is an element of $\left\{
1,2,\ldots,m\right\}  $, then $A_{\bullet,\sim v}=\left(  a_{i,\mathbf{d}%
_{v}\left(  j\right)  }\right)  _{1\leq i\leq n,\ 1\leq j\leq m-1}$.
\par
\textit{Observation 2:} If $n\in\mathbb{N}$ and $m\in\mathbb{N}$, if
$A=\left(  a_{i,j}\right)  _{1\leq i\leq n,\ 1\leq j\leq m}\in\mathbb{K}%
^{n\times m}$ is an $n\times m$-matrix, and if $w$ and $v$ are two elements of
$\left\{  1,2,\ldots,m\right\}  $ satisfying $w<v$, then $\operatorname*{cols}%
\nolimits_{1,2,\ldots,\widehat{w},\ldots,\widehat{v},\ldots,m}A=\left(
a_{i,\mathbf{d}_{w,v}\left(  j\right)  }\right)  _{1\leq i\leq n,\ 1\leq j\leq
m-2}$.
\par
\textit{Observation 3:} For every two integers $w$ and $v$ satisfying $w<v$,
we have
\[
\mathbf{d}_{w,v}\left(  j\right)  =\mathbf{d}_{v}\left(  \mathbf{d}_{w}\left(
j\right)  \right)  \ \ \ \ \ \ \ \ \ \ \text{for each }j\in\mathbb{Z}.
\]
\par
\textit{Proof of Observation 1:} Observation 1 was already proven in a
previous footnote (namely, the one where we gave a rigorous proof of
Proposition \ref{prop.unrows.basics} \textbf{(d)}).
\par
\textit{Proof of Observation 2:} Let $n\in\mathbb{N}$ and $m\in\mathbb{N}$.
Let $A=\left(  a_{i,j}\right)  _{1\leq i\leq n,\ 1\leq j\leq m}\in
\mathbb{K}^{n\times m}$ be an $n\times m$-matrix. Let $w$ and $v$ be two
elements of $\left\{  1,2,\ldots,m\right\}  $ satisfying $w<v$. Recalling how
$\mathbf{d}_{w,v}\left(  j\right)  $ is defined for every $j\in\mathbb{Z}$, we
find that%
\begin{align*}
&  \left(  \mathbf{d}_{w,v}\left(  1\right)  ,\mathbf{d}_{w,v}\left(
2\right)  ,\ldots,\mathbf{d}_{w,v}\left(  m-2\right)  \right) \\
&  =\left(  1,2,\ldots,w-1,w+1,w+2,\ldots,v-1,v+1,v+2,\ldots,m\right) \\
&  =\left(  1,2,\ldots,\widehat{w},\ldots,\widehat{v},\ldots,m\right)  .
\end{align*}
Now,%
\begin{align*}
&  \operatorname*{cols}\nolimits_{1,2,\ldots,\widehat{w},\ldots,\widehat{v}%
,\ldots,m}A\\
&  =\operatorname*{cols}\nolimits_{\mathbf{d}_{w,v}\left(  1\right)
,\mathbf{d}_{w,v}\left(  2\right)  ,\ldots,\mathbf{d}_{w,v}\left(  m-2\right)
}A\\
&  \ \ \ \ \ \ \ \ \ \ \left(  \text{since }\left(  1,2,\ldots,\widehat{w}%
,\ldots,\widehat{v},\ldots,m\right)  =\left(  \mathbf{d}_{w,v}\left(
1\right)  ,\mathbf{d}_{w,v}\left(  2\right)  ,\ldots,\mathbf{d}_{w,v}\left(
m-2\right)  \right)  \right) \\
&  =\left(  a_{i,\mathbf{d}_{w,v}\left(  y\right)  }\right)  _{1\leq i\leq
n,\ 1\leq y\leq m-2}\\
&  \ \ \ \ \ \ \ \ \ \ \left(
\begin{array}
[c]{c}%
\text{by the definition of }\operatorname*{cols}\nolimits_{\mathbf{d}%
_{w,v}\left(  1\right)  ,\mathbf{d}_{w,v}\left(  2\right)  ,\ldots
,\mathbf{d}_{w,v}\left(  m-2\right)  }A\\
\text{(since }A=\left(  a_{i,j}\right)  _{1\leq i\leq n,\ 1\leq j\leq
m}\text{)}%
\end{array}
\right) \\
&  =\left(  a_{i,\mathbf{d}_{w,v}\left(  j\right)  }\right)  _{1\leq i\leq
n,\ 1\leq j\leq m-2}\ \ \ \ \ \ \ \ \ \ \left(  \text{here, we have renamed
the index }\left(  i,y\right)  \text{ as }\left(  i,j\right)  \right)  .
\end{align*}
This proves Observation 2.
\par
\textit{Proof of Observation 3:} Let $w$ and $v$ be two integers satisfying
$w<v$. Then, $w\leq v-1$ (since $w$ and $v$ are integers). For every
$j\in\mathbb{Z}$, we can prove $\mathbf{d}_{w,v}\left(  j\right)
=\mathbf{d}_{v}\left(  \mathbf{d}_{w}\left(  j\right)  \right)  $ by an
explicit computation using the definitions of $\mathbf{d}_{w,v}\left(
j\right)  $, of $\mathbf{d}_{w}\left(  j\right)  $ and of $\mathbf{d}%
_{v}\left(  \mathbf{d}_{w}\left(  j\right)  \right)  $. (The cases when $j<w$,
when $w\leq j<v-1$ and when $v-1\leq j$ need to be treated separately; but
each of these three cases is completely straightforward. For instance, in the
case when $v-1\leq j$, the definition of $\mathbf{d}_{w,v}\left(  j\right)  $
yields%
\[
\mathbf{d}_{w,v}\left(  j\right)  =%
\begin{cases}
j, & \text{if }j<w;\\
j+1, & \text{if }w\leq j<v-1;\\
j+2, & \text{if }v-1\leq j
\end{cases}
=j+2\ \ \ \ \ \ \ \ \ \ \left(  \text{since }v-1\leq j\right)  ,
\]
whereas the definition of $\mathbf{d}_{w}\left(  j\right)  $ yields%
\[
\mathbf{d}_{w}\left(  j\right)  =%
\begin{cases}
j, & \text{if }j<w;\\
j+1, & \text{if }j\geq w
\end{cases}
=j+1\ \ \ \ \ \ \ \ \ \ \left(  \text{since }j\geq v-1\geq w\right)
\]
and therefore%
\begin{align*}
\mathbf{d}_{v}\left(  \mathbf{d}_{w}\left(  j\right)  \right)   &
=\mathbf{d}_{v}\left(  j+1\right)  =%
\begin{cases}
j+1, & \text{if }j+1<v;\\
\left(  j+1\right)  +1, & \text{if }j+1\geq v
\end{cases}
\ \ \ \ \ \ \ \ \ \ \left(  \text{by the definition of }\mathbf{d}_{v}\left(
j+1\right)  \right) \\
&  =\left(  j+1\right)  +1\ \ \ \ \ \ \ \ \ \ \left(  \text{since }j+1\geq
v\text{ (since }j\geq v-1\text{)}\right) \\
&  =j+2=\mathbf{d}_{w,v}\left(  j\right)  ;
\end{align*}
thus, $\mathbf{d}_{w,v}\left(  j\right)  =\mathbf{d}_{v}\left(  \mathbf{d}%
_{w}\left(  j\right)  \right)  $ is proven in this case.) This proves
Observation 3.
\par
Let us now prove Proposition \ref{prop.unrows.basics} \textbf{(h)}. Indeed,
let $v\in\left\{  1,2,\ldots,m\right\}  $ and $w\in\left\{  1,2,\ldots
,v-1\right\}  $. Then, $w\leq v-1$ (since $w\in\left\{  1,2,\ldots
,v-1\right\}  $) and thus $w\leq\underbrace{v}_{\leq m}-1\leq m-1$, so that
$v\in\left\{  1,2,\ldots,m-1\right\}  $.
\par
Write the matrix $A$ in the form $A=\left(  a_{i,j}\right)  _{1\leq i\leq
n,\ 1\leq j\leq m}$. Then, Observation 1 shows that $A_{\bullet,\sim
v}=\left(  a_{i,\mathbf{d}_{v}\left(  j\right)  }\right)  _{1\leq i\leq
n,\ 1\leq j\leq m-1}$. Hence, Observation 1 (applied to $m-1$, $A_{\bullet
,\sim v}$, $a_{i,\mathbf{d}_{v}\left(  j\right)  }$ and $w$ instead of $m$,
$A$, $a_{i,j}$ and $v$) yields%
\begin{align*}
\left(  A_{\bullet,\sim v}\right)  _{\bullet,\sim w}  &  =\left(
a_{i,\mathbf{d}_{v}\left(  \mathbf{d}_{w}\left(  j\right)  \right)  }\right)
_{1\leq i\leq n,\ 1\leq j\leq\left(  m-1\right)  -1}=\left(
\underbrace{a_{i,\mathbf{d}_{v}\left(  \mathbf{d}_{w}\left(  j\right)
\right)  }}_{\substack{=a_{i,\mathbf{d}_{w,v}\left(  j\right)  }\\\text{(since
Observation 3 yields}\\\mathbf{d}_{v}\left(  \mathbf{d}_{w}\left(  j\right)
\right)  =\mathbf{d}_{w,v}\left(  j\right)  \text{)}}}\right)  _{1\leq i\leq
n,\ 1\leq j\leq m-2}\\
&  =\left(  a_{i,\mathbf{d}_{w,v}\left(  j\right)  }\right)  _{1\leq i\leq
n,\ 1\leq j\leq m-2}=\operatorname*{cols}\nolimits_{1,2,\ldots,\widehat{w}%
,\ldots,\widehat{v},\ldots,m}A
\end{align*}
(by Observation 2). This proves Proposition \ref{prop.unrows.basics}
\textbf{(h)}.}

Similarly, we can prove Proposition \ref{prop.unrows.basics} \textbf{(i)}%
\footnote{We will not actually need Proposition \ref{prop.unrows.basics}
\textbf{(i)} in the following, so you can just as well skip this proof. The
same applies to Proposition \ref{prop.unrows.basics} \textbf{(k)}.}: Let
$v\in\left\{  1,2,\ldots,m\right\}  $ and $w\in\left\{  v,v+1,\ldots
,m-1\right\}  $. We need to show that $\left(  A_{\bullet,\sim v}\right)
_{\bullet,\sim w}=\operatorname*{cols}\nolimits_{1,2,\ldots,\widehat{v}%
,\ldots,\widehat{w+1},\ldots,m}A$. The matrix $\left(  A_{\bullet,\sim
v}\right)  _{\bullet,\sim w}$ is obtained from the matrix $A$ by first
removing the $v$-th column (thus obtaining an $n\times\left(  m-1\right)
$-matrix) and then removing the $w$-th column (from the resulting
$n\times\left(  m-1\right)  $-matrix). Since $w\geq v$ (because $w\in\left\{
v,v+1,\ldots,m-1\right\}  $), the first of these two removals (i.e., the
removal of the $v$-th column of $A$) has moved the $\left(  w+1\right)  $-th
column of $A$ one step to the left; therefore, the latter column has become
the $w$-th column after this first removal. The second removal then removes
this column. Therefore, the two successive removals could be replaced by a
simultaneous removal of both the $v$-th and the $\left(  w+1\right)  $-th
columns from the matrix $A$. In other words, the matrix obtained from the
matrix $A$ by first removing the $v$-th column and then removing the $w$-th
column is identical with the matrix obtained from $A$ by simultaneously
removing both the $v$-th and the $\left(  w+1\right)  $-th columns. But the
former matrix is $\left(  A_{\bullet,\sim v}\right)  _{\bullet,\sim w}$,
whereas the latter matrix is $\operatorname*{cols}\nolimits_{1,2,\ldots
,\widehat{v},\ldots,\widehat{w+1},\ldots,m}A$. Hence, the identity of these
two matrices rewrites as $\left(  A_{\bullet,\sim v}\right)  _{\bullet,\sim
w}=\operatorname*{cols}\nolimits_{1,2,\ldots,\widehat{v},\ldots,\widehat{w+1}%
,\ldots,m}A$. This proves Proposition \ref{prop.unrows.basics} \textbf{(i)}%
.\footnote{Again, a more rigorous proof can be given (if desired) similarly to
our rigorous proof of Proposition \ref{prop.unrows.basics} \textbf{(h)}.}

Proposition \ref{prop.unrows.basics} \textbf{(j)} is the analogue of
Proposition \ref{prop.unrows.basics} \textbf{(h)} for rows instead of columns
(with $v$ renamed as $u$); its proof is equally analogous.

Proposition \ref{prop.unrows.basics} \textbf{(k)} is the analogue of
Proposition \ref{prop.unrows.basics} \textbf{(i)} for rows instead of columns
(with $v$ renamed as $u$); its proof is equally analogous.

Let us now prove Proposition \ref{prop.unrows.basics} \textbf{(l)}: Let
$v\in\left\{  1,2,\ldots,n\right\}  $, $u\in\left\{  1,2,\ldots,n\right\}  $
and $q\in\left\{  1,2,\ldots,m\right\}  $ be such that $u<v$. We have $u<v$
and thus $u\leq v-1$ (since $u$ and $v$ are integers). Combined with $u\geq1$
(since $u\in\left\{  1,2,\ldots,n\right\}  $), this yields $u\in\left\{
1,2,\ldots,v-1\right\}  $.

Proposition \ref{prop.unrows.basics} \textbf{(c)} (applied to $v$ and $q$
instead of $u$ and $v$) yields $\left(  A_{\bullet,\sim q}\right)  _{\sim
v,\bullet}=\left(  A_{\sim v,\bullet}\right)  _{\bullet,\sim q}=A_{\sim v,\sim
q}$.

Now, $u\leq\underbrace{v}_{\substack{\leq n\\\text{(since }v\in\left\{
1,2,\ldots,n\right\}  \text{)}}}-1\leq n-1$. Combining this with $u\geq1$, we
obtain $u\in\left\{  1,2,\ldots,n-1\right\}  $. Hence, Proposition
\ref{prop.unrows.basics} \textbf{(c)} (applied to $n-1$, $A_{\sim v,\bullet}$
and $q$ instead of $n$, $A$ and $v$) yields $\left(  \left(  A_{\sim
v,\bullet}\right)  _{\bullet,\sim q}\right)  _{\sim u,\bullet}=\left(  \left(
A_{\sim v,\bullet}\right)  _{\sim u,\bullet}\right)  _{\bullet,\sim q}=\left(
A_{\sim v,\bullet}\right)  _{\sim u,\sim q}$. Hence,%
\begin{align*}
\left(  A_{\sim v,\bullet}\right)  _{\sim u,\sim q}  &  =\left(
\underbrace{\left(  A_{\sim v,\bullet}\right)  _{\bullet,\sim q}}_{=\left(
A_{\bullet,\sim q}\right)  _{\sim v,\bullet}}\right)  _{\sim u,\bullet
}=\left(  \left(  A_{\bullet,\sim q}\right)  _{\sim v,\bullet}\right)  _{\sim
u,\bullet}\\
&  =\operatorname*{rows}\nolimits_{1,2,\ldots,\widehat{u},\ldots
,\widehat{v},\ldots,n}\left(  A_{\bullet,\sim q}\right)
\end{align*}
(by Proposition \ref{prop.unrows.basics} \textbf{(j)}, applied to $m-1$,
$A_{\bullet,\sim q}$, $v$ and $u$ instead of $m$, $A$, $u$ and $w$). This
proves Proposition \ref{prop.unrows.basics} \textbf{(l)}.
\end{proof}

Thus, the proof of Proposition \ref{prop.unrows.basics} is complete. It
remains to prove Proposition \ref{prop.unrows.basics-I}. Let me give an
informal proof:

\begin{proof}
[Proof of Proposition \ref{prop.unrows.basics-I}.]We know that $\left(
I_{n}\right)  _{\bullet,u}$ is the $u$-th column of the identity matrix. All
entries of this column are $0$, except for the $u$-th, which is $1$. In other
words,%
\[
\left(  I_{n}\right)  _{\bullet,u}=\left(  \underbrace{0,0,\ldots
,0}_{u-1\text{ zeroes}},1,\underbrace{0,0,\ldots,0}_{n-u\text{ zeroes}%
}\right)  ^{T}.
\]
Now, $\left(  \left(  I_{n}\right)  _{\bullet,u}\right)  _{\sim v,\bullet}$ is
the result of removing the $v$-th row from this column, i.e., removing the
$v$-th entry from this column\footnote{Indeed, the $v$-th row is the same as
the $v$-th entry in this case, because each row of $\left(  I_{n}\right)
_{\bullet,u}$ consists of one entry only.}. This $v$-th row is below the $1$
in the $u$-th row (since $u<v$); therefore, the result of its removal is the
column vector
\[
\left(  \underbrace{0,0,\ldots,0}_{u-1\text{ zeroes}},1,\underbrace{0,0,\ldots
,0}_{n-u-1\text{ zeroes}}\right)  ^{T}.
\]
Hence, we have shown that%
\begin{align*}
\left(  \left(  I_{n}\right)  _{\bullet,u}\right)  _{\sim v,\bullet}  &
=\left(  \underbrace{0,0,\ldots,0}_{u-1\text{ zeroes}}%
,1,\underbrace{0,0,\ldots,0}_{n-u-1\text{ zeroes}}\right)  ^{T}\\
&  =\left(  \underbrace{0,0,\ldots,0}_{u-1\text{ zeroes}}%
,1,\underbrace{0,0,\ldots,0}_{\left(  n-1\right)  -u\text{ zeroes}}\right)
^{T}\\
&  =\left(  \text{the }u\text{-th column of the matrix }I_{n-1}\right)
=\left(  I_{n-1}\right)  _{\bullet,u}.
\end{align*}
This proves Proposition \ref{prop.unrows.basics-I}\footnote{Again, if this
proof was not rigorous enough for you, here is a more formal proof of
Proposition \ref{prop.unrows.basics-I}:
\par
For every $j\in\mathbb{Z}$ and $r\in\mathbb{Z}$, let $\mathbf{d}_{r}\left(
j\right)  $ be the integer $%
\begin{cases}
j, & \text{if }j<r;\\
j+1, & \text{if }j\geq r
\end{cases}
$. For any two objects $i$ and $j$, we define an element $\delta_{i,j}%
\in\mathbb{K}$ by $\delta_{i,j}=%
\begin{cases}
1, & \text{if }i=j;\\
0, & \text{if }i\neq j
\end{cases}
$.
\par
We make the following two observations:
\par
\textit{Observation 1:} If $n\in\mathbb{N}$ and $m\in\mathbb{N}$, if
$A=\left(  a_{i,j}\right)  _{1\leq i\leq n,\ 1\leq j\leq m}\in\mathbb{K}%
^{n\times m}$ is an $n\times m$-matrix, and if $u$ is an element of $\left\{
1,2,\ldots,n\right\}  $, then $A_{\sim u,\bullet}=\left(  a_{\mathbf{d}%
_{u}\left(  i\right)  ,j}\right)  _{1\leq i\leq n-1,\ 1\leq j\leq m}$.
\par
\textit{Observation 2:} If $u\in\mathbb{Z}$, $v\in\mathbb{Z}$ and
$i\in\mathbb{Z}$ satisfy $u<v$, then $\delta_{\mathbf{d}_{v}\left(  i\right)
,u}=\delta_{i,u}$.
\par
\textit{Proof of Observation 1:} Observation 1 is analogous to the Observation
1 from the footnote in our above proof of Proposition \ref{prop.unrows.basics}
\textbf{(h)}. (More precisely, it is the analogue of the latter observation
for rows instead of columns.)
\par
\textit{Proof of Observation 2:} Let $u\in\mathbb{Z}$, $v\in\mathbb{Z}$ and
$i\in\mathbb{Z}$ satisfy $u<v$. We are in one of the following two Cases:
\par
\textit{Case 1:} We have $i<v$.
\par
\textit{Case 2:} We have $i\geq v$.
\par
Let us first consider Case 1. In this case, we have $i<v$. Now, the definition
of $\mathbf{d}_{v}\left(  i\right)  $ yields $\mathbf{d}_{v}\left(  i\right)
=%
\begin{cases}
i, & \text{if }i<v;\\
i+1, & \text{if }i\geq v
\end{cases}
=i$ (since $i<v$) and thus $\delta_{\mathbf{d}_{v}\left(  i\right)  ,u}%
=\delta_{i,u}$. Hence, Observation 2 is proven in Case 1.
\par
Let us now consider Case 2. In this case, we have $i\geq v$. Now, the
definition of $\mathbf{d}_{v}\left(  i\right)  $ yields $\mathbf{d}_{v}\left(
i\right)  =%
\begin{cases}
i, & \text{if }i<v;\\
i+1, & \text{if }i\geq v
\end{cases}
=i+1$ (since $i\geq v$) and thus $\mathbf{d}_{v}\left(  i\right)  =i+1>i\geq
v>u$. Hence, $\mathbf{d}_{v}\left(  i\right)  \neq u$. Thus, $\delta
_{\mathbf{d}_{v}\left(  i\right)  ,u}=0$. Also, $i\geq v>u$ and thus $i\neq
u$; hence, $\delta_{i,u}=0$. Now, $\delta_{\mathbf{d}_{v}\left(  i\right)
,u}=0=\delta_{i,u}$. Hence, Observation 2 is proven in Case 2.
\par
Thus, Observation 2 is proven in both Cases 1 and 2; hence, Observation 2 is
always proven.
\par
Now, let $n$, $u$ and $v$ be as in Proposition \ref{prop.unrows.basics-I}.
From $u<v\leq n$, we obtain $u\in\left\{  1,2,\ldots,n-1\right\}  $, so that
$\left(  I_{n-1}\right)  _{\bullet,u}$ is well-defined. Also, $\left(  \left(
I_{n}\right)  _{\bullet,u}\right)  _{\sim v,\bullet}$ is well-defined (since
$u$ and $v$ belong to $\left\{  1,2,\ldots,n\right\}  $). Now, the definition
of $\left(  I_{n}\right)  _{\bullet,u}$ yields%
\begin{align*}
\left(  I_{n}\right)  _{\bullet,u}  &  =\left(  \text{the }u\text{-th column
of the matrix }I_{n}\right) \\
&  =\left(
\begin{array}
[c]{c}%
\delta_{1,u}\\
\delta_{2,u}\\
\vdots\\
\delta_{n,u}%
\end{array}
\right)  \ \ \ \ \ \ \ \ \ \ \left(  \text{since }I_{n}=\left(  \delta
_{i,j}\right)  _{1\leq i\leq n,\ 1\leq j\leq n}\right) \\
&  =\left(  \delta_{i,u}\right)  _{1\leq i\leq n,\ 1\leq j\leq1}.
\end{align*}
The same argument (applied to $n-1$ instead of $n$) yields $\left(
I_{n-1}\right)  _{\bullet,u}=\left(  \delta_{i,u}\right)  _{1\leq i\leq
n-1,\ 1\leq j\leq1}$.
\par
Recall that $\left(  I_{n}\right)  _{\bullet,u}=\left(  \delta_{i,u}\right)
_{1\leq i\leq n,\ 1\leq j\leq1}$. Hence, Observation 1 (applied to $1$,
$\left(  I_{n}\right)  _{\bullet,u}$, $\delta_{i,u}$ and $v$ instead of $m$,
$A$, $a_{i,j}$ and $u$) yields%
\[
\left(  \left(  I_{n}\right)  _{\bullet,u}\right)  _{\sim v,\bullet}=\left(
\underbrace{\delta_{\mathbf{d}_{v}\left(  i\right)  ,u}}_{\substack{=\delta
_{i,u}\\\text{(by Observation 2)}}}\right)  _{1\leq i\leq n-1,\ 1\leq j\leq
1}=\left(  \delta_{i,u}\right)  _{1\leq i\leq n-1,\ 1\leq j\leq1}.
\]
Comparing this with $\left(  I_{n-1}\right)  _{\bullet,u}=\left(  \delta
_{i,u}\right)  _{1\leq i\leq n-1,\ 1\leq j\leq1}$, we obtain $\left(  \left(
I_{n}\right)  _{\bullet,u}\right)  _{\sim v,\bullet}=\left(  I_{n-1}\right)
_{\bullet,u}$. This proves Proposition \ref{prop.unrows.basics-I}.}.
\end{proof}

We have now proven both Proposition \ref{prop.unrows.basics} and Proposition
\ref{prop.unrows.basics-I}. Thus, Exercise \ref{exe.unrows.basics}
\textbf{(a)} is solved.

\textbf{(b)} We need to derive Proposition \ref{prop.desnanot.12} and
Proposition \ref{prop.desnanot.1n} from Theorem \ref{thm.desnanot}.

\begin{proof}
[Proof of Proposition \ref{prop.desnanot.12} using Theorem \ref{thm.desnanot}%
.]We have%
\begin{equation}
\operatorname*{sub}\nolimits_{1,2,\ldots,\widehat{1},\ldots,\widehat{2}%
,\ldots,n}^{1,2,\ldots,\widehat{1},\ldots,\widehat{2},\ldots,n}A=\widetilde{A}
\label{pf.prop.desnanot.12.short.1}%
\end{equation}
\footnote{\textit{Proof of (\ref{pf.prop.desnanot.12.short.1}):} The $\left(
n-2\right)  $-tuple $\left(  1,2,\ldots,\widehat{1},\ldots,\widehat{2}%
,\ldots,n\right)  $ is obtained by removing the $1$-st and the $2$-nd entries
from the $n$-tuple $\left(  1,2,\ldots,n\right)  $. Thus,%
\[
\left(  1,2,\ldots,\widehat{1},\ldots,\widehat{2},\ldots,n\right)  =\left(
3,4,\ldots,n\right)  =\left(  1+2,2+2,\ldots,\left(  n-2\right)  +2\right)  .
\]
Hence,%
\begin{align*}
\operatorname*{sub}\nolimits_{1,2,\ldots,\widehat{1},\ldots,\widehat{2}%
,\ldots,n}^{1,2,\ldots,\widehat{1},\ldots,\widehat{2},\ldots,n}A  &
=\operatorname*{sub}\nolimits_{1+2,2+2,\ldots,\left(  n-2\right)
+2}^{1+2,2+2,\ldots,\left(  n-2\right)  +2}A=\left(  a_{x+2,y+2}\right)
_{1\leq x\leq n-2,\ 1\leq y\leq n-2}\\
&  \ \ \ \ \ \ \ \ \ \ \left(  \text{by the definition of }\operatorname*{sub}%
\nolimits_{1+2,2+2,\ldots,\left(  n-2\right)  +2}^{1+2,2+2,\ldots,\left(
n-2\right)  +2}A\right) \\
&  =\left(  a_{i+2,j+2}\right)  _{1\leq i\leq n-2,\ 1\leq j\leq n-2}\\
&  \ \ \ \ \ \ \ \ \ \ \left(  \text{here, we have renamed the index }\left(
x,y\right)  \text{ as }\left(  i,j\right)  \right) \\
&  =\widetilde{A}\ \ \ \ \ \ \ \ \ \ \left(  \text{since }\widetilde{A}%
=\left(  a_{i+2,j+2}\right)  _{1\leq i\leq n-2,\ 1\leq j\leq n-2}\right)  .
\end{align*}
Qed.}. Now, $1<2$. Moreover, $1$ and $2$ are elements of $\left\{
1,2,\ldots,n\right\}  $ (since $n\geq2$). Hence, Theorem \ref{thm.desnanot}
(applied to $p=1$, $q=2$, $u=1$ and $v=2$) yields%
\begin{align*}
&  \det A\cdot\det\left(  \operatorname*{sub}\nolimits_{1,2,\ldots
,\widehat{1},\ldots,\widehat{2},\ldots,n}^{1,2,\ldots,\widehat{1}%
,\ldots,\widehat{2},\ldots,n}A\right) \\
&  =\det\left(  A_{\sim1,\sim1}\right)  \cdot\det\left(  A_{\sim2,\sim
2}\right)  -\det\left(  A_{\sim1,\sim2}\right)  \cdot\det\left(  A_{\sim
2,\sim1}\right)  .
\end{align*}
In view of (\ref{pf.prop.desnanot.12.short.1}), this rewrites as
\begin{align*}
&  \det A\cdot\det\widetilde{A}\\
&  =\det\left(  A_{\sim1,\sim1}\right)  \cdot\det\left(  A_{\sim2,\sim
2}\right)  -\det\left(  A_{\sim1,\sim2}\right)  \cdot\det\left(  A_{\sim
2,\sim1}\right)  .
\end{align*}
This proves Proposition \ref{prop.desnanot.12}.
\end{proof}

\begin{proof}
[Proof of Proposition \ref{prop.desnanot.1n} using Theorem \ref{thm.desnanot}%
.]We have%
\begin{equation}
\operatorname*{sub}\nolimits_{1,2,\ldots,\widehat{1},\ldots,\widehat{n}%
,\ldots,n}^{1,2,\ldots,\widehat{1},\ldots,\widehat{n},\ldots,n}A=A^{\prime}
\label{pf.prop.desnanot.1n.short.1}%
\end{equation}
\footnote{\textit{Proof of (\ref{pf.prop.desnanot.1n.short.1}):} The $\left(
n-2\right)  $-tuple $\left(  1,2,\ldots,\widehat{1},\ldots,\widehat{n}%
,\ldots,n\right)  $ is obtained by removing the $1$-st and the $n$-th entries
from the $n$-tuple $\left(  1,2,\ldots,n\right)  $. Thus,%
\[
\left(  1,2,\ldots,\widehat{1},\ldots,\widehat{n},\ldots,n\right)  =\left(
2,3,\ldots,n-1\right)  =\left(  1+1,2+1,\ldots,\left(  n-2\right)  +1\right)
.
\]
Hence,%
\begin{align*}
\operatorname*{sub}\nolimits_{1,2,\ldots,\widehat{1},\ldots,\widehat{n}%
,\ldots,n}^{1,2,\ldots,\widehat{1},\ldots,\widehat{n},\ldots,n}A  &
=\operatorname*{sub}\nolimits_{1+1,2+1,\ldots,\left(  n-2\right)
+1}^{1+1,2+1,\ldots,\left(  n-2\right)  +1}A=\left(  a_{x+1,y+1}\right)
_{1\leq x\leq n-2,\ 1\leq y\leq n-2}\\
&  \ \ \ \ \ \ \ \ \ \ \left(  \text{by the definition of }\operatorname*{sub}%
\nolimits_{1+1,2+1,\ldots,\left(  n-2\right)  +1}^{1+1,2+1,\ldots,\left(
n-2\right)  +1}A\right) \\
&  =\left(  a_{i+1,j+1}\right)  _{1\leq i\leq n-2,\ 1\leq j\leq n-2}\\
&  \ \ \ \ \ \ \ \ \ \ \left(  \text{here, we have renamed the index }\left(
x,y\right)  \text{ as }\left(  i,j\right)  \right) \\
&  =A^{\prime}\ \ \ \ \ \ \ \ \ \ \left(  \text{since }A^{\prime}=\left(
a_{i+1,j+1}\right)  _{1\leq i\leq n-2,\ 1\leq j\leq n-2}\right)  .
\end{align*}
Qed.}. Now, $1<n$ (since $n\geq2$) and $1<n$. Moreover, $1$ and $n$ are
elements of $\left\{  1,2,\ldots,n\right\}  $ (since $n\geq2\geq1$). Hence,
Theorem \ref{thm.desnanot} (applied to $p=1$, $q=n$, $u=1$ and $v=n$) yields%
\begin{align*}
&  \det A\cdot\det\left(  \operatorname*{sub}\nolimits_{1,2,\ldots
,\widehat{1},\ldots,\widehat{n},\ldots,n}^{1,2,\ldots,\widehat{1}%
,\ldots,\widehat{n},\ldots,n}A\right) \\
&  =\det\left(  A_{\sim1,\sim1}\right)  \cdot\det\left(  A_{\sim n,\sim
n}\right)  -\det\left(  A_{\sim1,\sim n}\right)  \cdot\det\left(  A_{\sim
n,\sim1}\right)  .
\end{align*}
In view of (\ref{pf.prop.desnanot.1n.short.1}), this rewrites as
\begin{align*}
&  \det A\cdot\det\left(  A^{\prime}\right) \\
&  =\det\left(  A_{\sim1,\sim1}\right)  \cdot\det\left(  A_{\sim n,\sim
n}\right)  -\det\left(  A_{\sim1,\sim n}\right)  \cdot\det\left(  A_{\sim
n,\sim1}\right)  .
\end{align*}
This proves Proposition \ref{prop.desnanot.1n}.
\end{proof}

Now, both Proposition \ref{prop.desnanot.12} and Proposition
\ref{prop.desnanot.1n} are proven using Theorem \ref{thm.desnanot}. This
completes the solution of Exercise \ref{exe.unrows.basics} \textbf{(b)}.
\end{proof}
\end{vershort}

\begin{verlong}
In preparation for our solution to Exercise \ref{exe.unrows.basics}, let us
introduce the following notations:

\begin{definition}
\label{def.sol.unrows.d}\textbf{(a)} For every $j\in\mathbb{Z}$ and
$r\in\mathbb{Z}$, let $\mathbf{d}_{r}\left(  j\right)  $ be the integer $%
\begin{cases}
j, & \text{if }j<r;\\
j+1, & \text{if }j\geq r
\end{cases}
$.

\textbf{(b)} For every $j\in\mathbb{Z}$, $r\in\mathbb{Z}$ and $s\in\mathbb{Z}$
satisfying $r<s$, we let $\mathbf{d}_{r,s}\left(  j\right)  $ be the integer $%
\begin{cases}
j, & \text{if }j<r;\\
j+1, & \text{if }r\leq j<s-1;\\
j+2, & \text{if }s-1\leq j
\end{cases}
$. (This is well-defined, because $r\leq s-1$ (since $r<s$ and since $r$ and
$s$ are integers).)
\end{definition}

Using these notations, we can state three really basic facts:

\begin{proposition}
\label{prop.sol.unirows.r}Let $n\in\mathbb{N}$. Let $r\in\left\{
1,2,\ldots,n\right\}  $. Then,
\[
\left(  \mathbf{d}_{r}\left(  1\right)  ,\mathbf{d}_{r}\left(  2\right)
,\ldots,\mathbf{d}_{r}\left(  n-1\right)  \right)  =\left(  1,2,\ldots
,\widehat{r},\ldots,n\right)  .
\]

\end{proposition}

\begin{proof}
[Proof of Proposition \ref{prop.sol.unirows.r}.]We have $\mathbf{d}_{r}\left(
j\right)  =%
\begin{cases}
j, & \text{if }j<r;\\
j+1, & \text{if }j\geq r
\end{cases}
$ for each $j\in\left\{  1,2,\ldots,n-1\right\}  $ (by the definition of
$\mathbf{d}_{r}\left(  j\right)  $). From this, we conclude that
\begin{equation}
\left(  \mathbf{d}_{r}\left(  1\right)  ,\mathbf{d}_{r}\left(  2\right)
,\ldots,\mathbf{d}_{r}\left(  n-1\right)  \right)  =\left(  1,2,\ldots
,r-1,r+1,r+2,\ldots,n\right)  \label{pf.prop.sol.unirows.r.1}%
\end{equation}
\footnote{Here is this argument in more detail:
\par
If $\alpha$ and $\beta$ are two finite lists of integers, then we define a new
finite list $\alpha\ast\beta$ of integers by setting%
\[
\alpha\ast\beta=\left(  \alpha_{1},\alpha_{2},\ldots,\alpha_{a},\beta
_{1},\beta_{2},\ldots,\beta_{b}\right)  ,
\]
where $\alpha$ is written in the form $\left(  \alpha_{1},\alpha_{2}%
,\ldots,\alpha_{a}\right)  $, and where $\beta$ is written in the form
$\left(  \beta_{1},\beta_{2},\ldots,\beta_{b}\right)  $. This new list
$\alpha\ast\beta$ is called the \textit{concatenation} of the lists $\alpha$
and $\beta$.
\par
For every $j\in\left\{  1,2,\ldots,r-1\right\}  $, we have%
\begin{align*}
\mathbf{d}_{r}\left(  j\right)   &  =%
\begin{cases}
j, & \text{if }j<r;\\
j+1, & \text{if }j\geq r
\end{cases}
\ \ \ \ \ \ \ \ \ \ \left(  \text{by the definition of }\mathbf{d}_{r}\left(
j\right)  \right) \\
&  =j\ \ \ \ \ \ \ \ \ \ \left(  \text{since }j<r\text{ (because }j\in\left\{
1,2,\ldots,r-1\right\}  \text{)}\right)  .
\end{align*}
In other words, $\left(  \mathbf{d}_{r}\left(  1\right)  ,\mathbf{d}%
_{r}\left(  2\right)  ,\ldots,\mathbf{d}_{r}\left(  r-1\right)  \right)
=\left(  1,2,\ldots,r-1\right)  $.
\par
On the other hand, for every $j\in\left\{  r,r+1,\ldots,n-1\right\}  $, we
have%
\begin{align*}
\mathbf{d}_{r}\left(  j\right)   &  =%
\begin{cases}
j, & \text{if }j<r;\\
j+1, & \text{if }j\geq r
\end{cases}
\ \ \ \ \ \ \ \ \ \ \left(  \text{by the definition of }\mathbf{d}_{r}\left(
j\right)  \right) \\
&  =j+1\ \ \ \ \ \ \ \ \ \ \left(  \text{since }j\geq r\text{ (because }%
j\in\left\{  r,r+1,\ldots,n-1\right\}  \text{)}\right)  .
\end{align*}
In other words, $\left(  \mathbf{d}_{r}\left(  r\right)  ,\mathbf{d}%
_{r}\left(  r+1\right)  ,\ldots,\mathbf{d}_{r}\left(  n-1\right)  \right)
=\left(  r+1,\left(  r+1\right)  +1,\ldots,\left(  n-1\right)  +1\right)  $.
Now,%
\begin{align*}
&  \left(  \mathbf{d}_{r}\left(  1\right)  ,\mathbf{d}_{r}\left(  2\right)
,\ldots,\mathbf{d}_{r}\left(  n-1\right)  \right) \\
&  =\left(  \mathbf{d}_{r}\left(  1\right)  ,\mathbf{d}_{r}\left(  2\right)
,\ldots,\mathbf{d}_{r}\left(  r-1\right)  ,\mathbf{d}_{r}\left(  r\right)
,\mathbf{d}_{r}\left(  r+1\right)  ,\ldots,\mathbf{d}_{r}\left(  n-1\right)
\right) \\
&  \ \ \ \ \ \ \ \ \ \ \left(  \text{since }r\in\left\{  1,2,\ldots,n\right\}
\right) \\
&  =\underbrace{\left(  \mathbf{d}_{r}\left(  1\right)  ,\mathbf{d}_{r}\left(
2\right)  ,\ldots,\mathbf{d}_{r}\left(  r-1\right)  \right)  }_{=\left(
1,2,\ldots,r-1\right)  }\ast\underbrace{\left(  \mathbf{d}_{r}\left(
r\right)  ,\mathbf{d}_{r}\left(  r+1\right)  ,\ldots,\mathbf{d}_{r}\left(
n-1\right)  \right)  }_{\substack{=\left(  r+1,\left(  r+1\right)
+1,\ldots,\left(  n-1\right)  +1\right)  \\=\left(  r+1,r+2,\ldots,n\right)
}}\\
&  =\left(  1,2,\ldots,r-1\right)  \ast\left(  r+1,r+2,\ldots,n\right)
=\left(  1,2,\ldots,r-1,r+1,r+2,\ldots,n\right)
\end{align*}
(since $r\in\left\{  1,2,\ldots,n\right\}  $). This proves
(\ref{pf.prop.sol.unirows.r.1}).}. Now, the definition of $\left(
1,2,\ldots,\widehat{r},\ldots,n\right)  $ yields%
\begin{align*}
\left(  1,2,\ldots,\widehat{r},\ldots,n\right)   &  =\left(  1,2,\ldots
,r-1,r+1,r+2,\ldots,n\right) \\
&  =\left(  \mathbf{d}_{r}\left(  1\right)  ,\mathbf{d}_{r}\left(  2\right)
,\ldots,\mathbf{d}_{r}\left(  n-1\right)  \right)  \ \ \ \ \ \ \ \ \ \ \left(
\text{by (\ref{pf.prop.sol.unirows.r.1})}\right)  .
\end{align*}
This proves Proposition \ref{prop.sol.unirows.r}.
\end{proof}

\begin{proposition}
\label{prop.sol.unirows.d}Let $n\in\mathbb{N}$ and $m\in\mathbb{N}$. Let
$A=\left(  a_{i,j}\right)  _{1\leq i\leq n,\ 1\leq j\leq m}\in\mathbb{K}%
^{n\times m}$ be an $n\times m$-matrix.

\textbf{(a)} For every $u\in\left\{  1,2,\ldots,n\right\}  $, we have $A_{\sim
u,\bullet}=\left(  a_{\mathbf{d}_{u}\left(  i\right)  ,j}\right)  _{1\leq
i\leq n-1,\ 1\leq j\leq m}$.

\textbf{(b)} For every $v\in\left\{  1,2,\ldots,m\right\}  $, we have
$A_{\bullet,\sim v}=\left(  a_{i,\mathbf{d}_{v}\left(  j\right)  }\right)
_{1\leq i\leq n,\ 1\leq j\leq m-1}$.
\end{proposition}

\begin{proof}
[Proof of Proposition \ref{prop.sol.unirows.d}.]\textbf{(a)} Let $u\in\left\{
1,2,\ldots,n\right\}  $. Then, Proposition \ref{prop.sol.unirows.r} (applied
to $r=u$) yields%
\[
\left(  \mathbf{d}_{u}\left(  1\right)  ,\mathbf{d}_{u}\left(  2\right)
,\ldots,\mathbf{d}_{u}\left(  n-1\right)  \right)  =\left(  1,2,\ldots
,\widehat{u},\ldots,n\right)  .
\]

The definition of $A_{\sim u,\bullet}$ yields
\begin{align*}
A_{\sim u,\bullet}  &  =\operatorname*{rows}\nolimits_{1,2,\ldots
,\widehat{u},\ldots,n}A=\operatorname*{rows}\nolimits_{\mathbf{d}_{u}\left(
1\right)  ,\mathbf{d}_{u}\left(  2\right)  ,\ldots,\mathbf{d}_{u}\left(
n-1\right)  }A\\
&  \ \ \ \ \ \ \ \ \ \ \left(  \text{since }\left(  1,2,\ldots,\widehat{u}%
,\ldots,n\right)  =\left(  \mathbf{d}_{u}\left(  1\right)  ,\mathbf{d}%
_{u}\left(  2\right)  ,\ldots,\mathbf{d}_{u}\left(  n-1\right)  \right)
\right) \\
&  =\left(  a_{\mathbf{d}_{u}\left(  x\right)  ,j}\right)  _{1\leq x\leq
n-1,\ 1\leq j\leq m}\\
&  \ \ \ \ \ \ \ \ \ \ \left(
\begin{array}
[c]{c}%
\text{by the definition of }\operatorname*{rows}\nolimits_{\mathbf{d}%
_{u}\left(  1\right)  ,\mathbf{d}_{u}\left(  2\right)  ,\ldots,\mathbf{d}%
_{u}\left(  n-1\right)  }A\\
\text{(since }A=\left(  a_{i,j}\right)  _{1\leq i\leq n,\ 1\leq j\leq
m}\text{)}%
\end{array}
\right) \\
&  =\left(  a_{\mathbf{d}_{u}\left(  i\right)  ,j}\right)  _{1\leq i\leq
n-1,\ 1\leq j\leq m}\\
&  \ \ \ \ \ \ \ \ \ \ \left(  \text{here, we have renamed the index }\left(
x,j\right)  \text{ as }\left(  i,j\right)  \right)  .
\end{align*}
This proves Proposition \ref{prop.sol.unirows.d} \textbf{(a)}.

\textbf{(b)} Let $v\in\left\{  1,2,\ldots,m\right\}  $. Then, Proposition
\ref{prop.sol.unirows.r} (applied to $v$ and $m$ instead of $r$ and $n$)
yields%
\[
\left(  \mathbf{d}_{v}\left(  1\right)  ,\mathbf{d}_{v}\left(  2\right)
,\ldots,\mathbf{d}_{v}\left(  m-1\right)  \right)  =\left(  1,2,\ldots
,\widehat{v},\ldots,m\right)  .
\]

The definition of $A_{\bullet,\sim v}$ yields
\begin{align*}
A_{\bullet,\sim v}  &  =\operatorname*{cols}\nolimits_{1,2,\ldots
,\widehat{v},\ldots,m}A=\operatorname*{cols}\nolimits_{\mathbf{d}_{v}\left(
1\right)  ,\mathbf{d}_{v}\left(  2\right)  ,\ldots,\mathbf{d}_{v}\left(
m-1\right)  }A\\
&  \ \ \ \ \ \ \ \ \ \ \left(  \text{since }\left(  1,2,\ldots,\widehat{v}%
,\ldots,m\right)  =\left(  \mathbf{d}_{v}\left(  1\right)  ,\mathbf{d}%
_{v}\left(  2\right)  ,\ldots,\mathbf{d}_{v}\left(  m-1\right)  \right)
\right) \\
&  =\left(  a_{i,\mathbf{d}_{v}\left(  y\right)  }\right)  _{1\leq i\leq
n,\ 1\leq y\leq m-1}\\
&  \ \ \ \ \ \ \ \ \ \ \left(
\begin{array}
[c]{c}%
\text{by the definition of }\operatorname*{cols}\nolimits_{\mathbf{d}%
_{v}\left(  1\right)  ,\mathbf{d}_{v}\left(  2\right)  ,\ldots,\mathbf{d}%
_{v}\left(  m-1\right)  }A\\
\text{(since }A=\left(  a_{i,j}\right)  _{1\leq i\leq n,\ 1\leq j\leq
m}\text{)}%
\end{array}
\right) \\
&  =\left(  a_{i,\mathbf{d}_{v}\left(  j\right)  }\right)  _{1\leq i\leq
n,\ 1\leq j\leq m-1}\\
&  \ \ \ \ \ \ \ \ \ \ \left(  \text{here, we have renamed the index }\left(
i,y\right)  \text{ as }\left(  i,j\right)  \right)  .
\end{align*}
This proves Proposition \ref{prop.sol.unirows.d} \textbf{(b)}.
\end{proof}

\begin{proposition}
\label{prop.sol.unrows.2}Let $n\in\mathbb{N}$. Let $r$ and $s$ be two elements
of $\left\{  1,2,\ldots,n\right\}  $ such that $r<s$.

\textbf{(a)} Then,
\[
\left(  1,2,\ldots,\widehat{r},\ldots,\widehat{s},\ldots,n\right)  =\left(
\mathbf{d}_{r,s}\left(  1\right)  ,\mathbf{d}_{r,s}\left(  2\right)
,\ldots,\mathbf{d}_{r,s}\left(  n-2\right)  \right)  .
\]

\textbf{(b)} We have
\[
\mathbf{d}_{r,s}\left(  j\right)  =\mathbf{d}_{s}\left(  \mathbf{d}_{r}\left(
j\right)  \right)  \ \ \ \ \ \ \ \ \ \ \text{for every }j\in\mathbb{Z}.
\]

\textbf{(c)} We have
\[
\mathbf{d}_{r,s}\left(  j\right)  =\mathbf{d}_{r}\left(  \mathbf{d}%
_{s-1}\left(  j\right)  \right)  \ \ \ \ \ \ \ \ \ \ \text{for every }%
j\in\mathbb{Z}.
\]

\end{proposition}

\begin{proof}
[Proof of Proposition \ref{prop.sol.unrows.2}.]Recall that $r<s$. Thus, $r\leq
s-1$ (since $r$ and $s$ are integers). Hence, $\underbrace{r}_{\leq s-1}%
-1\leq\left(  s-1\right)  -1=s-2$.

\textbf{(a)} If $\alpha$ and $\beta$ are two finite lists of integers, then we
define a new finite list $\alpha\ast\beta$ of integers by setting%
\[
\alpha\ast\beta=\left(  \alpha_{1},\alpha_{2},\ldots,\alpha_{a},\beta
_{1},\beta_{2},\ldots,\beta_{b}\right)  ,
\]
where $\alpha$ is written in the form $\left(  \alpha_{1},\alpha_{2}%
,\ldots,\alpha_{a}\right)  $, and where $\beta$ is written in the form
$\left(  \beta_{1},\beta_{2},\ldots,\beta_{b}\right)  $. This new list
$\alpha\ast\beta$ is called the \textit{concatenation} of the lists $\alpha$
and $\beta$. If $\alpha$, $\beta$ and $\gamma$ are three finite lists of
integers, then the two lists $\left(  \alpha\ast\beta\right)  \ast\gamma$ and
$\alpha\ast\left(  \beta\ast\gamma\right)  $ are identical. In other words,
concatenation of finite lists is an associative operation. Thus, we are able
to write $\alpha\ast\beta\ast\gamma$ for any of the two identical lists
$\left(  \alpha\ast\beta\right)  \ast\gamma$ and $\alpha\ast\left(  \beta
\ast\gamma\right)  $ whenever $\alpha$, $\beta$ and $\gamma$ are three finite
lists of integers.

Now, the definition of $\left(  1,2,\ldots,\widehat{r},\ldots,\widehat{s}%
,\ldots,n\right)  $ yields%
\begin{align}
&  \left(  1,2,\ldots,\widehat{r},\ldots,\widehat{s},\ldots,n\right)
\nonumber\\
&  =\left(  \underbrace{1,2,\ldots,r-1}_{\substack{\text{all integers}%
\\\text{from }1\text{ to }r-1}},\underbrace{r+1,r+2,\ldots,s-1}%
_{\substack{\text{all integers}\\\text{from }r+1\text{ to }s-1}%
},\underbrace{s+1,s+2,\ldots,n}_{\substack{\text{all integers}\\\text{from
}s+1\text{ to }n}}\right) \nonumber\\
&  =\underbrace{\left(  \underbrace{1,2,\ldots,r-1}_{\substack{\text{all
integers}\\\text{from }1\text{ to }r-1}},\underbrace{r+1,r+2,\ldots
,s-1}_{\substack{\text{all integers}\\\text{from }r+1\text{ to }s-1}}\right)
}_{=\left(  1,2,\ldots,r-1\right)  \ast\left(  r+1,r+2,\ldots,s-1\right)
}\ast\left(  s+1,s+2,\ldots,n\right) \nonumber\\
&  =\left(  1,2,\ldots,r-1\right)  \ast\left(  r+1,r+2,\ldots,s-1\right)
\ast\left(  s+1,s+2,\ldots,n\right)  . \label{pf.prop.sol.unrows.2.a.5}%
\end{align}
Now,
\begin{equation}
\mathbf{d}_{r,s}\left(  j\right)  =j\ \ \ \ \ \ \ \ \ \ \text{for every }%
j\in\left\{  1,2,\ldots,r-1\right\}  \label{pf.prop.sol.unrows.2.a.7a}%
\end{equation}
\footnote{\textit{Proof of (\ref{pf.prop.sol.unrows.2.a.7a}):} Let
$j\in\left\{  1,2,\ldots,r-1\right\}  $. Then, $j\in\mathbb{Z}$ and $j<r$
(since $j\in\left\{  1,2,\ldots,r-1\right\}  $). The definition of
$\mathbf{d}_{r,s}\left(  j\right)  $ yields $\mathbf{d}_{r,s}\left(  j\right)
=%
\begin{cases}
j, & \text{if }j<r;\\
j+1, & \text{if }r\leq j<s-1;\\
j+2, & \text{if }s-1\leq j
\end{cases}
=j$ (since $j<r$). This proves (\ref{pf.prop.sol.unrows.2.a.7a}).}. In other
words,%
\begin{equation}
\left(  \mathbf{d}_{r,s}\left(  1\right)  ,\mathbf{d}_{r,s}\left(  2\right)
,\ldots,\mathbf{d}_{r,s}\left(  r-1\right)  \right)  =\left(  1,2,\ldots
,r-1\right)  . \label{pf.prop.sol.unrows.2.a.7b}%
\end{equation}

Also,%
\begin{equation}
\mathbf{d}_{r,s}\left(  j\right)  =j+1\ \ \ \ \ \ \ \ \ \ \text{for every
}j\in\left\{  r,r+1,\ldots,s-2\right\}  \label{pf.prop.sol.unrows.2.a.8a}%
\end{equation}
\footnote{\textit{Proof of (\ref{pf.prop.sol.unrows.2.a.8a}):} Let
$j\in\left\{  r,r+1,\ldots,s-2\right\}  $. Then, $j\in\mathbb{Z}$. Also,
$j\in\left\{  r,r+1,\ldots,s-2\right\}  $, so that $r\leq j\leq s-2$. From
$j\leq s-2<s-1$ and $r\leq j$, we obtain $r\leq j<s-1$. The definition of
$\mathbf{d}_{r,s}\left(  j\right)  $ yields $\mathbf{d}_{r,s}\left(  j\right)
=%
\begin{cases}
j, & \text{if }j<r;\\
j+1, & \text{if }r\leq j<s-1;\\
j+2, & \text{if }s-1\leq j
\end{cases}
=j+1$ (since $r\leq j<s-1$). This proves (\ref{pf.prop.sol.unrows.2.a.8a}).}.
In other words,%
\begin{align}
\left(  \mathbf{d}_{r,s}\left(  r\right)  ,\mathbf{d}_{r,s}\left(  r+1\right)
,\ldots,\mathbf{d}_{r,s}\left(  s-2\right)  \right)   &  =\left(  r+1,\left(
r+1\right)  +1,\ldots,\left(  s-2\right)  +1\right) \nonumber\\
&  =\left(  r+1,r+2,\ldots,s-1\right)  . \label{pf.prop.sol.unrows.2.a.8b}%
\end{align}

Finally,%
\begin{equation}
\mathbf{d}_{r,s}\left(  j\right)  =j+2\ \ \ \ \ \ \ \ \ \ \text{for every
}j\in\left\{  s-1,s,\ldots,n-2\right\}  \label{pf.prop.sol.unrows.2.a.9a}%
\end{equation}
\footnote{\textit{Proof of (\ref{pf.prop.sol.unrows.2.a.9a}):} Let
$j\in\left\{  s-1,s,\ldots,n-2\right\}  $. Then, $j\in\mathbb{Z}$. Also,
$j\in\left\{  s-1,s,\ldots,n-2\right\}  $, so that $j\geq s-1$. Thus, $s-1\leq
j$. The definition of $\mathbf{d}_{r,s}\left(  j\right)  $ yields
$\mathbf{d}_{r,s}\left(  j\right)  =%
\begin{cases}
j, & \text{if }j<r;\\
j+1, & \text{if }r\leq j<s-1;\\
j+2, & \text{if }s-1\leq j
\end{cases}
=j+2$ (since $s-1\leq j$). This proves (\ref{pf.prop.sol.unrows.2.a.9a}).}. In
other words,%
\begin{align}
\left(  \mathbf{d}_{r,s}\left(  s-1\right)  ,\mathbf{d}_{r,s}\left(  s\right)
,\ldots,\mathbf{d}_{r,s}\left(  n-2\right)  \right)   &  =\left(  \left(
s-1\right)  +2,s+2,\ldots,\left(  n-2\right)  +2\right) \nonumber\\
&  =\left(  s+1,s+2,\ldots,n\right)  . \label{pf.prop.sol.unrows.2.a.9b}%
\end{align}

Now, recall that $r-1\leq s-2$. Hence,%
\begin{align}
&  \left(  \mathbf{d}_{r,s}\left(  1\right)  ,\mathbf{d}_{r,s}\left(
2\right)  ,\ldots,\mathbf{d}_{r,s}\left(  s-2\right)  \right) \nonumber\\
&  =\underbrace{\left(  \mathbf{d}_{r,s}\left(  1\right)  ,\mathbf{d}%
_{r,s}\left(  2\right)  ,\ldots,\mathbf{d}_{r,s}\left(  r-1\right)  \right)
}_{\substack{=\left(  1,2,\ldots,r-1\right)  \\\text{(by
(\ref{pf.prop.sol.unrows.2.a.7b}))}}}\ast\underbrace{\left(  \mathbf{d}%
_{r,s}\left(  r\right)  ,\mathbf{d}_{r,s}\left(  r+1\right)  ,\ldots
,\mathbf{d}_{r,s}\left(  s-2\right)  \right)  }_{\substack{=\left(
r+1,r+2,\ldots,s-1\right)  \\\text{(by (\ref{pf.prop.sol.unrows.2.a.8b}))}%
}}\nonumber\\
&  =\left(  1,2,\ldots,r-1\right)  \ast\left(  r+1,r+2,\ldots,s-1\right)  .
\label{pf.prop.sol.unrows.2.a.10b}%
\end{align}
Now,%
\begin{align*}
&  \left(  \mathbf{d}_{r,s}\left(  1\right)  ,\mathbf{d}_{r,s}\left(
2\right)  ,\ldots,\mathbf{d}_{r,s}\left(  n-2\right)  \right) \\
&  =\underbrace{\left(  \mathbf{d}_{r,s}\left(  1\right)  ,\mathbf{d}%
_{r,s}\left(  2\right)  ,\ldots,\mathbf{d}_{r,s}\left(  s-2\right)  \right)
}_{\substack{=\left(  1,2,\ldots,r-1\right)  \ast\left(  r+1,r+2,\ldots
,s-1\right)  \\\text{(by (\ref{pf.prop.sol.unrows.2.a.10b}))}}}\ast
\underbrace{\left(  \mathbf{d}_{r,s}\left(  s-1\right)  ,\mathbf{d}%
_{r,s}\left(  s\right)  ,\ldots,\mathbf{d}_{r,s}\left(  n-2\right)  \right)
}_{\substack{=\left(  s+1,s+2,\ldots,n\right)  \\\text{(by
(\ref{pf.prop.sol.unrows.2.a.9b}))}}}\\
&  =\left(  1,2,\ldots,r-1\right)  \ast\left(  r+1,r+2,\ldots,s-1\right)
\ast\left(  s+1,s+2,\ldots,n\right) \\
&  =\left(  1,2,\ldots,\widehat{r},\ldots,\widehat{s},\ldots,n\right)
\ \ \ \ \ \ \ \ \ \ \left(  \text{by (\ref{pf.prop.sol.unrows.2.a.5})}\right)
.
\end{align*}
This proves Proposition \ref{prop.sol.unrows.2} \textbf{(a)}.

\textbf{(b)} Let $j\in\mathbb{Z}$. We want to prove that $\mathbf{d}%
_{r,s}\left(  j\right)  =\mathbf{d}_{s}\left(  \mathbf{d}_{r}\left(  j\right)
\right)  $.

We have either $j<s-1$ or $j\geq s-1$. In other words, we must be in one of
the following two Cases:

\textit{Case 1:} We have $j<s-1$.

\textit{Case 2:} We have $j\geq s-1$.

Let us first consider Case 1. In this case, we have $j<s-1$. Thus, $j<s-1<s$.
We have either $j<r$ or $j\geq r$. In other words, we must be in one of the
following two Subcases:

\textit{Subcase 1.1:} We have $j<r$.

\textit{Subcase 1.2:} We have $j\geq r$.

Let us first consider Subcase 1.1. In this subcase, we have $j<r$. Now, the
definition of $\mathbf{d}_{r}\left(  j\right)  $ yields $\mathbf{d}_{r}\left(
j\right)  =%
\begin{cases}
j, & \text{if }j<r;\\
j+1, & \text{if }j\geq r
\end{cases}
=j$ (since $j<r$). But $j<r<s$. Now,
\begin{align*}
\mathbf{d}_{s}\left(  \underbrace{\mathbf{d}_{r}\left(  j\right)  }%
_{=j}\right)   &  =\mathbf{d}_{s}\left(  j\right)  =%
\begin{cases}
j, & \text{if }j<s;\\
j+1, & \text{if }j\geq s
\end{cases}
\ \ \ \ \ \ \ \ \ \ \left(  \text{by the definition of }\mathbf{d}_{s}\left(
j\right)  \right) \\
&  =j\ \ \ \ \ \ \ \ \ \ \left(  \text{since }j<s\right)  .
\end{align*}
Comparing this with%
\begin{align*}
\mathbf{d}_{r,s}\left(  j\right)   &  =%
\begin{cases}
j, & \text{if }j<r;\\
j+1, & \text{if }r\leq j<s-1;\\
j+2, & \text{if }s-1\leq j
\end{cases}
\ \ \ \ \ \ \ \ \ \ \left(  \text{by the definition of }\mathbf{d}%
_{r,s}\left(  j\right)  \right) \\
&  =j\ \ \ \ \ \ \ \ \ \ \left(  \text{since }j<r\right)  ,
\end{align*}
we obtain $\mathbf{d}_{r,s}\left(  j\right)  =\mathbf{d}_{s}\left(
\mathbf{d}_{r}\left(  j\right)  \right)  $. Thus, $\mathbf{d}_{r,s}\left(
j\right)  =\mathbf{d}_{s}\left(  \mathbf{d}_{r}\left(  j\right)  \right)  $ is
proven in Subcase 1.1.

Let us now consider Subcase 1.2. In this subcase, we have $j\geq r$. In other
words, $r\leq j$. Now, the definition of $\mathbf{d}_{r}\left(  j\right)  $
yields $\mathbf{d}_{r}\left(  j\right)  =%
\begin{cases}
j, & \text{if }j<r;\\
j+1, & \text{if }j\geq r
\end{cases}
=j+1$ (since $j\geq r$). But $j<s-1$, so that $j+1<s$. Now,
\begin{align*}
\mathbf{d}_{s}\left(  \underbrace{\mathbf{d}_{r}\left(  j\right)  }%
_{=j+1}\right)   &  =\mathbf{d}_{s}\left(  j+1\right) \\
&  =%
\begin{cases}
j+1, & \text{if }j+1<s;\\
\left(  j+1\right)  +1, & \text{if }j+1\geq s
\end{cases}
\ \ \ \ \ \ \ \ \ \ \left(  \text{by the definition of }\mathbf{d}_{s}\left(
j+1\right)  \right) \\
&  =j+1\ \ \ \ \ \ \ \ \ \ \left(  \text{since }j+1<s\right)  .
\end{align*}
Comparing this with%
\begin{align*}
\mathbf{d}_{r,s}\left(  j\right)   &  =%
\begin{cases}
j, & \text{if }j<r;\\
j+1, & \text{if }r\leq j<s-1;\\
j+2, & \text{if }s-1\leq j
\end{cases}
\ \ \ \ \ \ \ \ \ \ \left(  \text{by the definition of }\mathbf{d}%
_{r,s}\left(  j\right)  \right) \\
&  =j\ \ \ \ \ \ \ \ \ \ \left(  \text{since }r\leq j<s-1\text{ (because
}r\leq j\text{ and }j<s-1\text{)}\right)  ,
\end{align*}
we obtain $\mathbf{d}_{r,s}\left(  j\right)  =\mathbf{d}_{s}\left(
\mathbf{d}_{r}\left(  j\right)  \right)  $. Thus, $\mathbf{d}_{r,s}\left(
j\right)  =\mathbf{d}_{s}\left(  \mathbf{d}_{r}\left(  j\right)  \right)  $ is
proven in Subcase 1.2.

We have now proven $\mathbf{d}_{r,s}\left(  j\right)  =\mathbf{d}_{s}\left(
\mathbf{d}_{r}\left(  j\right)  \right)  $ in each of the two Subcases 1.1 and
1.2. Since these two Subcases cover all of Case 1, we thus conclude that
$\mathbf{d}_{r,s}\left(  j\right)  =\mathbf{d}_{s}\left(  \mathbf{d}%
_{r}\left(  j\right)  \right)  $ holds in Case 1.

Let us now consider Case 2. In this case, we have $j\geq s-1$. Now, $r\leq
s-1\leq j$ (since $j\geq s-1$) and thus $j\geq r$. Now, the definition of
$\mathbf{d}_{r}\left(  j\right)  $ yields $\mathbf{d}_{r}\left(  j\right)  =%
\begin{cases}
j, & \text{if }j<r;\\
j+1, & \text{if }j\geq r
\end{cases}
=j+1$ (since $j\geq r$). But $j\geq s-1$, so that $j+1\geq s$. Now,
\begin{align*}
\mathbf{d}_{s}\left(  \underbrace{\mathbf{d}_{r}\left(  j\right)  }%
_{=j+1}\right)   &  =\mathbf{d}_{s}\left(  j+1\right) \\
&  =%
\begin{cases}
j+1, & \text{if }j+1<s;\\
\left(  j+1\right)  +1, & \text{if }j+1\geq s
\end{cases}
\ \ \ \ \ \ \ \ \ \ \left(  \text{by the definition of }\mathbf{d}_{s}\left(
j+1\right)  \right) \\
&  =\left(  j+1\right)  +1\ \ \ \ \ \ \ \ \ \ \left(  \text{since }j+1\geq
s\right) \\
&  =j+2.
\end{align*}
Comparing this with%
\begin{align*}
\mathbf{d}_{r,s}\left(  j\right)   &  =%
\begin{cases}
j, & \text{if }j<r;\\
j+1, & \text{if }r\leq j<s-1;\\
j+2, & \text{if }s-1\leq j
\end{cases}
\ \ \ \ \ \ \ \ \ \ \left(  \text{by the definition of }\mathbf{d}%
_{r,s}\left(  j\right)  \right) \\
&  =j+2\ \ \ \ \ \ \ \ \ \ \left(  \text{since }s-1\leq j\right)  ,
\end{align*}
we obtain $\mathbf{d}_{r,s}\left(  j\right)  =\mathbf{d}_{s}\left(
\mathbf{d}_{r}\left(  j\right)  \right)  $. Thus, $\mathbf{d}_{r,s}\left(
j\right)  =\mathbf{d}_{s}\left(  \mathbf{d}_{r}\left(  j\right)  \right)  $ is
proven in Case 2.

We have now proven $\mathbf{d}_{r,s}\left(  j\right)  =\mathbf{d}_{s}\left(
\mathbf{d}_{r}\left(  j\right)  \right)  $ in each of the two Cases 1 and 2.
Since these two Cases cover all possibilities, we thus conclude that
$\mathbf{d}_{r,s}\left(  j\right)  =\mathbf{d}_{s}\left(  \mathbf{d}%
_{r}\left(  j\right)  \right)  $ always holds. This completes the proof of
Proposition \ref{prop.sol.unrows.2} \textbf{(b)}.

\textbf{(c)} Let $j\in\mathbb{Z}$. We want to prove that $\mathbf{d}%
_{r,s}\left(  j\right)  =\mathbf{d}_{r}\left(  \mathbf{d}_{s-1}\left(
j\right)  \right)  $.

We have either $j<s-1$ or $j\geq s-1$. In other words, we must be in one of
the following two Cases:

\textit{Case 1:} We have $j<s-1$.

\textit{Case 2:} We have $j\geq s-1$.

Let us first consider Case 1. In this case, we have $j<s-1$. Thus, $j<s-1<s$.
Now, the definition of $\mathbf{d}_{s-1}\left(  j\right)  $ yields
$\mathbf{d}_{s-1}\left(  j\right)  =%
\begin{cases}
j, & \text{if }j<s-1;\\
j+1, & \text{if }j\geq s-1
\end{cases}
=j$ (since $j<s-1$).

We have either $j<r$ or $j\geq r$. In other words, we must be in one of the
following two Subcases:

\textit{Subcase 1.1:} We have $j<r$.

\textit{Subcase 1.2:} We have $j\geq r$.

Let us first consider Subcase 1.1. In this subcase, we have $j<r$. Now,%
\begin{align*}
\mathbf{d}_{r}\left(  \underbrace{\mathbf{d}_{s-1}\left(  j\right)  }%
_{=j}\right)   &  =\mathbf{d}_{r}\left(  j\right)  =%
\begin{cases}
j, & \text{if }j<r;\\
j+1, & \text{if }j\geq r
\end{cases}
\ \ \ \ \ \ \ \ \ \ \left(  \text{by the definition of }\mathbf{d}_{r}\left(
j\right)  \right) \\
&  =j\ \ \ \ \ \ \ \ \ \ \left(  \text{since }j<r\right)  .
\end{align*}
Comparing this with%
\begin{align*}
\mathbf{d}_{r,s}\left(  j\right)   &  =%
\begin{cases}
j, & \text{if }j<r;\\
j+1, & \text{if }r\leq j<s-1;\\
j+2, & \text{if }s-1\leq j
\end{cases}
\ \ \ \ \ \ \ \ \ \ \left(  \text{by the definition of }\mathbf{d}%
_{r,s}\left(  j\right)  \right) \\
&  =j\ \ \ \ \ \ \ \ \ \ \left(  \text{since }j<r\right)  ,
\end{align*}
we obtain $\mathbf{d}_{r,s}\left(  j\right)  =\mathbf{d}_{r}\left(
\mathbf{d}_{s-1}\left(  j\right)  \right)  $. Thus, $\mathbf{d}_{r,s}\left(
j\right)  =\mathbf{d}_{r}\left(  \mathbf{d}_{s-1}\left(  j\right)  \right)  $
is proven in Subcase 1.1.

Let us now consider Subcase 1.2. In this subcase, we have $j\geq r$. In other
words, $r\leq j$. Now,%
\begin{align*}
\mathbf{d}_{r}\left(  \underbrace{\mathbf{d}_{s-1}\left(  j\right)  }%
_{=j}\right)   &  =\mathbf{d}_{r}\left(  j\right)  =%
\begin{cases}
j, & \text{if }j<r;\\
j+1, & \text{if }j\geq r
\end{cases}
\ \ \ \ \ \ \ \ \ \ \left(  \text{by the definition of }\mathbf{d}_{r}\left(
j\right)  \right) \\
&  =j+1\ \ \ \ \ \ \ \ \ \ \left(  \text{since }j\geq r\right)  .
\end{align*}
Comparing this with%
\begin{align*}
\mathbf{d}_{r,s}\left(  j\right)   &  =%
\begin{cases}
j, & \text{if }j<r;\\
j+1, & \text{if }r\leq j<s-1;\\
j+2, & \text{if }s-1\leq j
\end{cases}
\ \ \ \ \ \ \ \ \ \ \left(  \text{by the definition of }\mathbf{d}%
_{r,s}\left(  j\right)  \right) \\
&  =j\ \ \ \ \ \ \ \ \ \ \left(  \text{since }r\leq j<s-1\text{ (because
}r\leq j\text{ and }j<s-1\text{)}\right)  ,
\end{align*}
we obtain $\mathbf{d}_{r,s}\left(  j\right)  =\mathbf{d}_{r}\left(
\mathbf{d}_{s-1}\left(  j\right)  \right)  $. Thus, $\mathbf{d}_{r,s}\left(
j\right)  =\mathbf{d}_{r}\left(  \mathbf{d}_{s-1}\left(  j\right)  \right)  $
is proven in Subcase 1.2.

We have now proven $\mathbf{d}_{r,s}\left(  j\right)  =\mathbf{d}_{r}\left(
\mathbf{d}_{s-1}\left(  j\right)  \right)  $ in each of the two Subcases 1.1
and 1.2. Since these two Subcases cover all of Case 1, we thus conclude that
$\mathbf{d}_{r,s}\left(  j\right)  =\mathbf{d}_{r}\left(  \mathbf{d}%
_{s-1}\left(  j\right)  \right)  $ holds in Case 1.

Let us now consider Case 2. In this case, we have $j\geq s-1$. In other words,
$s-1\leq j$. Now, the definition of $\mathbf{d}_{s-1}\left(  j\right)  $
yields $\mathbf{d}_{s-1}\left(  j\right)  =%
\begin{cases}
j, & \text{if }j<s-1;\\
j+1, & \text{if }j\geq s-1
\end{cases}
=j+1$ (since $j\geq s-1$). Also, $j\geq s-1$, so that $j+1\geq s>r$ and thus
$j+1\geq r$. Now,
\begin{align*}
\mathbf{d}_{r}\left(  \underbrace{\mathbf{d}_{s-1}\left(  j\right)  }%
_{=j+1}\right)   &  =\mathbf{d}_{r}\left(  j+1\right) \\
&  =%
\begin{cases}
j+1, & \text{if }j+1<r;\\
\left(  j+1\right)  +1, & \text{if }j+1\geq r
\end{cases}
\ \ \ \ \ \ \ \ \ \ \left(  \text{by the definition of }\mathbf{d}_{r}\left(
j+1\right)  \right) \\
&  =\left(  j+1\right)  +1\ \ \ \ \ \ \ \ \ \ \left(  \text{since }j+1\geq
r\right) \\
&  =j+2.
\end{align*}
Comparing this with%
\begin{align*}
\mathbf{d}_{r,s}\left(  j\right)   &  =%
\begin{cases}
j, & \text{if }j<r;\\
j+1, & \text{if }r\leq j<s-1;\\
j+2, & \text{if }s-1\leq j
\end{cases}
\ \ \ \ \ \ \ \ \ \ \left(  \text{by the definition of }\mathbf{d}%
_{r,s}\left(  j\right)  \right) \\
&  =j+2\ \ \ \ \ \ \ \ \ \ \left(  \text{since }s-1\leq j\right)  ,
\end{align*}
we obtain $\mathbf{d}_{r,s}\left(  j\right)  =\mathbf{d}_{r}\left(
\mathbf{d}_{s-1}\left(  j\right)  \right)  $. Thus, $\mathbf{d}_{r,s}\left(
j\right)  =\mathbf{d}_{r}\left(  \mathbf{d}_{s-1}\left(  j\right)  \right)  $
is proven in Case 2.

We have now proven $\mathbf{d}_{r,s}\left(  j\right)  =\mathbf{d}_{r}\left(
\mathbf{d}_{s-1}\left(  j\right)  \right)  $ in each of the two Cases 1 and 2.
Since these two Cases cover all possibilities, we thus conclude that
$\mathbf{d}_{r,s}\left(  j\right)  =\mathbf{d}_{r}\left(  \mathbf{d}%
_{s-1}\left(  j\right)  \right)  $ always holds. This completes the proof of
Proposition \ref{prop.sol.unrows.2} \textbf{(c)}.
\end{proof}

Now, we can step to the proof of Proposition \ref{prop.unrows.basics}:

\begin{proof}
[Proof of Proposition \ref{prop.unrows.basics}.]Write the matrix
$A\in\mathbb{K}^{n\times m}$ in the form $A=\left(  a_{i,j}\right)  _{1\leq
i\leq n,\ 1\leq j\leq m}$.

\textbf{(a)} Let $u\in\left\{  1,2,\ldots,n\right\}  $. The definition of
$A_{u,\bullet}$ shows that $A_{u,\bullet}$ is the $u$-th row of the matrix
$A$. In other words, $A_{u,\bullet}=\left(  \text{the }u\text{-th row of the
matrix }A\right)  $.

On the other hand, the definition of $\operatorname*{rows}\nolimits_{u}A$
yields $\operatorname*{rows}\nolimits_{u}A=\left(  a_{u,j}\right)  _{1\leq
x\leq1,\ 1\leq j\leq m}$. Now,%
\begin{align*}
A_{u,\bullet}  &  =\left(  \text{the }u\text{-th row of the matrix }A\right)
=\left(  a_{u,1},a_{u,2},\ldots,a_{u,m}\right) \\
&  \ \ \ \ \ \ \ \ \ \ \left(  \text{since }A=\left(  a_{i,j}\right)  _{1\leq
i\leq n,\ 1\leq j\leq m}\right) \\
&  =\left(  a_{u,j}\right)  _{1\leq x\leq1,\ 1\leq j\leq m}%
=\operatorname*{rows}\nolimits_{u}A\ \ \ \ \ \ \ \ \ \ \left(  \text{since
}\operatorname*{rows}\nolimits_{u}A=\left(  a_{u,j}\right)  _{1\leq
x\leq1,\ 1\leq j\leq m}\right)  .
\end{align*}
This proves Proposition \ref{prop.unrows.basics} \textbf{(a)}.

\textbf{(b)} Let $v\in\left\{  1,2,\ldots,m\right\}  $. The definition of
$A_{\bullet,v}$ shows that $A_{\bullet,v}$ is the $v$-th column of the matrix
$A$. In other words, $A_{\bullet,v}=\left(  \text{the }v\text{-th column of
the matrix }A\right)  $.

On the other hand, the definition of $\operatorname*{cols}\nolimits_{v}A$
yields $\operatorname*{cols}\nolimits_{v}A=\left(  a_{i,v}\right)  _{1\leq
i\leq n,\ 1\leq y\leq1}$. Now,%
\begin{align*}
A_{\bullet,v}  &  =\left(  \text{the }v\text{-th column of the matrix
}A\right)  =\left(
\begin{array}
[c]{c}%
a_{1,v}\\
a_{2,v}\\
\vdots\\
a_{n,v}%
\end{array}
\right) \\
&  \ \ \ \ \ \ \ \ \ \ \left(  \text{since }A=\left(  a_{i,j}\right)  _{1\leq
i\leq n,\ 1\leq j\leq m}\right) \\
&  =\left(  a_{i,v}\right)  _{1\leq i\leq n,\ 1\leq y\leq1}%
=\operatorname*{cols}\nolimits_{v}A\ \ \ \ \ \ \ \ \ \ \left(  \text{since
}\operatorname*{cols}\nolimits_{v}A=\left(  a_{i,v}\right)  _{1\leq i\leq
n,\ 1\leq y\leq1}\right)  .
\end{align*}
This proves Proposition \ref{prop.unrows.basics} \textbf{(b)}.

\textbf{(c)} Let $u\in\left\{  1,2,\ldots,n\right\}  $ and $v\in\left\{
1,2,\ldots,m\right\}  $. Then, Proposition \ref{prop.submatrix.easy}
\textbf{(d)} (applied to $n-1$, $\left(  1,2,\ldots,\widehat{u},\ldots
,n\right)  $, $m-1$ and $\left(  1,2,\ldots,\widehat{v},\ldots,m\right)  $
instead of $u$, $\left(  i_{1},i_{2},\ldots,i_{u}\right)  $, $v$ and $\left(
j_{1},j_{2},\ldots,j_{v}\right)  $) yields%
\begin{align}
\operatorname*{sub}\nolimits_{1,2,\ldots,\widehat{u},\ldots,n}^{1,2,\ldots
,\widehat{v},\ldots,m}A  &  =\operatorname*{rows}\nolimits_{1,2,\ldots
,\widehat{u},\ldots,n}\left(  \operatorname*{cols}\nolimits_{1,2,\ldots
,\widehat{v},\ldots,m}A\right) \label{pf.prop.unrows.basics.c.1}\\
&  =\operatorname*{cols}\nolimits_{1,2,\ldots,\widehat{v},\ldots,m}\left(
\operatorname*{rows}\nolimits_{1,2,\ldots,\widehat{u},\ldots,n}A\right)  .
\label{pf.prop.unrows.basics.c.2}%
\end{align}
But the definition of $A_{\sim u,\sim v}$ yields $A_{\sim u,\sim
v}=\operatorname*{sub}\nolimits_{1,2,\ldots,\widehat{u},\ldots,n}%
^{1,2,\ldots,\widehat{v},\ldots,m}A$.

Now, $A_{\sim u,\bullet}=\operatorname*{rows}\nolimits_{1,2,\ldots
,\widehat{u},\ldots,n}A$ (by the definition of $A_{\sim u,\bullet}$), and%
\begin{align}
\left(  A_{\sim u,\bullet}\right)  _{\bullet,\sim v}  &  =\operatorname*{cols}%
\nolimits_{1,2,\ldots,\widehat{v},\ldots,m}\left(  \underbrace{A_{\sim
u,\bullet}}_{=\operatorname*{rows}\nolimits_{1,2,\ldots,\widehat{u},\ldots
,n}A}\right)  \ \ \ \ \ \ \ \ \ \ \left(  \text{by the definition of }\left(
A_{\sim u,\bullet}\right)  _{\bullet,\sim v}\right) \nonumber\\
&  =\operatorname*{cols}\nolimits_{1,2,\ldots,\widehat{v},\ldots,m}\left(
\operatorname*{rows}\nolimits_{1,2,\ldots,\widehat{u},\ldots,n}A\right)
=\operatorname*{sub}\nolimits_{1,2,\ldots,\widehat{u},\ldots,n}^{1,2,\ldots
,\widehat{v},\ldots,m}A\ \ \ \ \ \ \ \ \ \ \left(  \text{by
(\ref{pf.prop.unrows.basics.c.2})}\right) \nonumber\\
&  =A_{\sim u,\sim v}\ \ \ \ \ \ \ \ \ \ \left(  \text{since }A_{\sim u,\sim
v}=\operatorname*{sub}\nolimits_{1,2,\ldots,\widehat{u},\ldots,n}%
^{1,2,\ldots,\widehat{v},\ldots,m}A\right)  .
\label{pf.prop.unrows.basics.c.5}%
\end{align}
Also, $A_{\bullet,\sim v}=\operatorname*{cols}\nolimits_{1,2,\ldots
,\widehat{v},\ldots,m}A$ (by the definition of $A_{\bullet,\sim v}$), and%
\begin{align*}
\left(  A_{\bullet,\sim v}\right)  _{\sim u,\bullet}  &  =\operatorname*{rows}%
\nolimits_{1,2,\ldots,\widehat{u},\ldots,n}\left(  \underbrace{A_{\bullet,\sim
v}}_{=\operatorname*{cols}\nolimits_{1,2,\ldots,\widehat{v},\ldots,m}%
A}\right)  \ \ \ \ \ \ \ \ \ \ \left(  \text{by the definition of }\left(
A_{\bullet,\sim v}\right)  _{\sim u,\bullet}\right) \\
&  =\operatorname*{rows}\nolimits_{1,2,\ldots,\widehat{u},\ldots,n}\left(
\operatorname*{cols}\nolimits_{1,2,\ldots,\widehat{v},\ldots,m}A\right)
=\operatorname*{sub}\nolimits_{1,2,\ldots,\widehat{u},\ldots,n}^{1,2,\ldots
,\widehat{v},\ldots,m}A\ \ \ \ \ \ \ \ \ \ \left(  \text{by
(\ref{pf.prop.unrows.basics.c.1})}\right) \\
&  =A_{\sim u,\sim v}\ \ \ \ \ \ \ \ \ \ \left(  \text{since }A_{\sim u,\sim
v}=\operatorname*{sub}\nolimits_{1,2,\ldots,\widehat{u},\ldots,n}%
^{1,2,\ldots,\widehat{v},\ldots,m}A\right)  .
\end{align*}
Combining this with (\ref{pf.prop.unrows.basics.c.5}), we obtain $\left(
A_{\bullet,\sim v}\right)  _{\sim u,\bullet}=\left(  A_{\sim u,\bullet
}\right)  _{\bullet,\sim v}=A_{\sim u,\sim v}$. This proves Proposition
\ref{prop.unrows.basics} \textbf{(c)}.

\textbf{(d)} Let $v\in\left\{  1,2,\ldots,m\right\}  $ and $w\in\left\{
1,2,\ldots,v-1\right\}  $. From $w\in\left\{  1,2,\ldots,v-1\right\}  $, we
obtain $w\geq1$ and $w\leq v-1$. Thus, $w\leq v-1<v$, so that $w\leq v\leq m$
(since $v\in\left\{  1,2,\ldots,m\right\}  $). Also, $w\in\left\{
1,2,\ldots,v-1\right\}  \subseteq\left\{  1,2,\ldots,m-1\right\}  $ (since
$\underbrace{v}_{\leq m}-1\leq m-1$). From $w\geq1$ and $w\leq m$, we obtain
$w\in\left\{  1,2,\ldots,m\right\}  $. The definition of $\mathbf{d}%
_{v}\left(  w\right)  $ yields $\mathbf{d}_{v}\left(  w\right)  =%
\begin{cases}
w, & \text{if }w<v;\\
w+1, & \text{if }w\geq v
\end{cases}
=w$ (since $w<v$).

But $A_{\bullet,\sim v}=\left(  a_{i,\mathbf{d}_{v}\left(  j\right)  }\right)
_{1\leq i\leq n,\ 1\leq j\leq m-1}$ (by Proposition \ref{prop.sol.unirows.d}
\textbf{(b)}). But $\left(  A_{\bullet,\sim v}\right)  _{\bullet,w}$ is the
$w$-th column of the matrix $A_{\bullet,\sim v}$ (by the definition of
$\left(  A_{\bullet,\sim v}\right)  _{\bullet,w}$). Thus,%
\begin{align*}
\left(  A_{\bullet,\sim v}\right)  _{\bullet,w}  &  =\left(  \text{the
}w\text{-th column of the matrix }A_{\bullet,\sim v}\right) \\
&  =\left(
\begin{array}
[c]{c}%
a_{1,\mathbf{d}_{v}\left(  w\right)  }\\
a_{2,\mathbf{d}_{v}\left(  w\right)  }\\
\vdots\\
a_{n,\mathbf{d}_{v}\left(  w\right)  }%
\end{array}
\right)  \ \ \ \ \ \ \ \ \ \ \left(  \text{since }A_{\bullet,\sim v}=\left(
a_{i,\mathbf{d}_{v}\left(  j\right)  }\right)  _{1\leq i\leq n,\ 1\leq j\leq
m-1}\right) \\
&  =\left(
\begin{array}
[c]{c}%
a_{1,w}\\
a_{2,w}\\
\vdots\\
a_{n,w}%
\end{array}
\right)  \ \ \ \ \ \ \ \ \ \ \left(  \text{since }\mathbf{d}_{v}\left(
w\right)  =w\right)  .
\end{align*}
Comparing this with%
\begin{align*}
A_{\bullet,w}  &  =\left(  \text{the }w\text{-th column of the matrix
}A\right) \\
&  \ \ \ \ \ \ \ \ \ \ \left(
\begin{array}
[c]{c}%
\text{since }A_{\bullet,w}\text{ is the }w\text{-th column of the matrix }A\\
\text{(by the definition of }A_{\bullet,w}\text{)}%
\end{array}
\right) \\
&  =\left(
\begin{array}
[c]{c}%
a_{1,w}\\
a_{2,w}\\
\vdots\\
a_{n,w}%
\end{array}
\right)  \ \ \ \ \ \ \ \ \ \ \left(  \text{since }A=\left(  a_{i,j}\right)
_{1\leq i\leq n,\ 1\leq j\leq m}\right)  ,
\end{align*}
we obtain $\left(  A_{\bullet,\sim v}\right)  _{\bullet,w}=A_{\bullet,w}$.
This proves Proposition \ref{prop.unrows.basics} \textbf{(d)}.

\textbf{(e)} Let $v\in\left\{  1,2,\ldots,m\right\}  $ and $w\in\left\{
v,v+1,\ldots,m-1\right\}  $. From $w\in\left\{  v,v+1,\ldots,m-1\right\}  $,
we obtain $w\geq v$ and $w\leq m-1$. Now, $w+1>w\geq v$, so that $w+1\geq
v\geq1$ (since $v\in\left\{  1,2,\ldots,m\right\}  $). Combining this with
$w+1\leq m$ (since $w\leq m-1$), we obtain $w+1\in\left\{  1,2,\ldots
,m\right\}  $. Also, $w\in\left\{  v,v+1,\ldots,m-1\right\}  \subseteq\left\{
1,2,\ldots,m-1\right\}  $ (since $v\geq1$).

The definition of $\mathbf{d}_{v}\left(  w\right)  $ yields $\mathbf{d}%
_{v}\left(  w\right)  =%
\begin{cases}
w, & \text{if }w<v;\\
w+1, & \text{if }w\geq v
\end{cases}
=w+1$ (since $w\geq v$).

But $A_{\bullet,\sim v}=\left(  a_{i,\mathbf{d}_{v}\left(  j\right)  }\right)
_{1\leq i\leq n,\ 1\leq j\leq m-1}$ (by Proposition \ref{prop.sol.unirows.d}
\textbf{(b)}). But $\left(  A_{\bullet,\sim v}\right)  _{\bullet,w}$ is the
$w$-th column of the matrix $A_{\bullet,\sim v}$ (by the definition of
$\left(  A_{\bullet,\sim v}\right)  _{\bullet,w}$). Thus,%
\begin{align*}
\left(  A_{\bullet,\sim v}\right)  _{\bullet,w}  &  =\left(  \text{the
}w\text{-th column of the matrix }A_{\bullet,\sim v}\right) \\
&  =\left(
\begin{array}
[c]{c}%
a_{1,\mathbf{d}_{v}\left(  w\right)  }\\
a_{2,\mathbf{d}_{v}\left(  w\right)  }\\
\vdots\\
a_{n,\mathbf{d}_{v}\left(  w\right)  }%
\end{array}
\right)  \ \ \ \ \ \ \ \ \ \ \left(  \text{since }A_{\bullet,\sim v}=\left(
a_{i,\mathbf{d}_{v}\left(  j\right)  }\right)  _{1\leq i\leq n,\ 1\leq j\leq
m-1}\right) \\
&  =\left(
\begin{array}
[c]{c}%
a_{1,w+1}\\
a_{2,w+1}\\
\vdots\\
a_{n,w+1}%
\end{array}
\right)  \ \ \ \ \ \ \ \ \ \ \left(  \text{since }\mathbf{d}_{v}\left(
w\right)  =w+1\right)  .
\end{align*}
Comparing this with%
\begin{align*}
A_{\bullet,w+1}  &  =\left(  \text{the }\left(  w+1\right)  \text{-th column
of the matrix }A\right) \\
&  \ \ \ \ \ \ \ \ \ \ \left(
\begin{array}
[c]{c}%
\text{since }A_{\bullet,w+1}\text{ is the }\left(  w+1\right)  \text{-th
column of the matrix }A\\
\text{(by the definition of }A_{\bullet,w+1}\text{)}%
\end{array}
\right) \\
&  =\left(
\begin{array}
[c]{c}%
a_{1,w+1}\\
a_{2,w+1}\\
\vdots\\
a_{n,w+1}%
\end{array}
\right)  \ \ \ \ \ \ \ \ \ \ \left(  \text{since }A=\left(  a_{i,j}\right)
_{1\leq i\leq n,\ 1\leq j\leq m}\right)  ,
\end{align*}
we obtain $\left(  A_{\bullet,\sim v}\right)  _{\bullet,w}=A_{\bullet,w+1}$.
This proves Proposition \ref{prop.unrows.basics} \textbf{(e)}.

\textbf{(f)} Let $u\in\left\{  1,2,\ldots,n\right\}  $ and $w\in\left\{
1,2,\ldots,u-1\right\}  $. From $w\in\left\{  1,2,\ldots,u-1\right\}  $, we
obtain $w\geq1$ and $w\leq u-1$. Thus, $w\leq u-1<u$, so that $w\leq u\leq n$
(since $u\in\left\{  1,2,\ldots,n\right\}  $). From $w\geq1$ and $w\leq n$, we
obtain $w\in\left\{  1,2,\ldots,n\right\}  $. We have $w\in\left\{
1,2,\ldots,u-1\right\}  \subseteq\left\{  1,2,\ldots,n-1\right\}  $ (since
$\underbrace{u}_{\leq n}-1\leq n-1$). The definition of $\mathbf{d}_{u}\left(
w\right)  $ yields $\mathbf{d}_{u}\left(  w\right)  =%
\begin{cases}
w, & \text{if }w<u;\\
w+1, & \text{if }w\geq u
\end{cases}
=w$ (since $w<u$).

But $A_{\sim u,\bullet}=\left(  a_{\mathbf{d}_{u}\left(  i\right)  ,j}\right)
_{1\leq i\leq n-1,\ 1\leq j\leq m}$ (by Proposition \ref{prop.sol.unirows.d}
\textbf{(a)}). But $\left(  A_{\sim u,\bullet}\right)  _{w,\bullet}$ is the
$w$-th row of the matrix $A_{\sim u,\bullet}$ (by the definition of $\left(
A_{\sim u,\bullet}\right)  _{w,\bullet}$). Thus,%
\begin{align*}
\left(  A_{\sim u,\bullet}\right)  _{w,\bullet}  &  =\left(  \text{the
}w\text{-th row of the matrix }A_{\sim u,\bullet}\right) \\
&  =\left(  a_{\mathbf{d}_{u}\left(  w\right)  ,1},a_{\mathbf{d}_{u}\left(
w\right)  ,2},\ldots,a_{\mathbf{d}_{u}\left(  w\right)  ,m}\right) \\
&  \ \ \ \ \ \ \ \ \ \ \left(  \text{since }A_{\sim u,\bullet}=\left(
a_{\mathbf{d}_{u}\left(  i\right)  ,j}\right)  _{1\leq i\leq n-1,\ 1\leq j\leq
m}\right) \\
&  =\left(  a_{w,1},a_{w,2},\ldots,a_{w,m}\right)  \ \ \ \ \ \ \ \ \ \ \left(
\text{since }\mathbf{d}_{u}\left(  w\right)  =w\right)  .
\end{align*}
Comparing this with%
\begin{align*}
A_{w,\bullet}  &  =\left(  \text{the }w\text{-th row of the matrix }A\right)
\\
&  \ \ \ \ \ \ \ \ \ \ \left(
\begin{array}
[c]{c}%
\text{since }A_{w,\bullet}\text{ is the }w\text{-th row of the matrix }A\\
\text{(by the definition of }A_{w,\bullet}\text{)}%
\end{array}
\right) \\
&  =\left(  a_{w,1},a_{w,2},\ldots,a_{w,m}\right)  \ \ \ \ \ \ \ \ \ \ \left(
\text{since }A=\left(  a_{i,j}\right)  _{1\leq i\leq n,\ 1\leq j\leq
m}\right)  ,
\end{align*}
we obtain $\left(  A_{\sim u,\bullet}\right)  _{w,\bullet}=A_{w,\bullet}$.
This proves Proposition \ref{prop.unrows.basics} \textbf{(f)}.

\textbf{(g)} Let $u\in\left\{  1,2,\ldots,n\right\}  $ and $w\in\left\{
u,u+1,\ldots,n-1\right\}  $. From $w\in\left\{  u,u+1,\ldots,n-1\right\}  $,
we obtain $w\geq u$ and $w\leq n-1$. Now, $w+1>w\geq u$, so that $w+1\geq
u\geq1$ (since $u\in\left\{  1,2,\ldots,n\right\}  $). Combining this with
$w+1\leq n$ (since $w\leq n-1$), we obtain $w+1\in\left\{  1,2,\ldots
,n\right\}  $. We have $w\in\left\{  u,u+1,\ldots,n-1\right\}  \subseteq
\left\{  1,2,\ldots,n-1\right\}  $ (since $u\geq1$).

The definition of $\mathbf{d}_{u}\left(  w\right)  $ yields $\mathbf{d}%
_{u}\left(  w\right)  =%
\begin{cases}
w, & \text{if }w<u;\\
w+1, & \text{if }w\geq u
\end{cases}
=w+1$ (since $w\geq u$).

But $A_{\sim u,\bullet}=\left(  a_{\mathbf{d}_{u}\left(  i\right)  ,j}\right)
_{1\leq i\leq n-1,\ 1\leq j\leq m}$ (by Proposition \ref{prop.sol.unirows.d}
\textbf{(a)}). But $\left(  A_{\sim u,\bullet}\right)  _{w,\bullet}$ is the
$w$-th row of the matrix $A_{\sim u,\bullet}$ (by the definition of $\left(
A_{\sim u,\bullet}\right)  _{w,\bullet}$). Thus,%
\begin{align*}
\left(  A_{\sim u,\bullet}\right)  _{w,\bullet}  &  =\left(  \text{the
}w\text{-th row of the matrix }A_{\sim u,\bullet}\right) \\
&  =\left(  a_{\mathbf{d}_{u}\left(  w\right)  ,1},a_{\mathbf{d}_{u}\left(
w\right)  ,2},\ldots,a_{\mathbf{d}_{u}\left(  w\right)  ,m}\right) \\
&  \ \ \ \ \ \ \ \ \ \ \left(  \text{since }A_{\sim u,\bullet}=\left(
a_{\mathbf{d}_{u}\left(  i\right)  ,j}\right)  _{1\leq i\leq n-1,\ 1\leq j\leq
m}\right) \\
&  =\left(  a_{w+1,1},a_{w+1,2},\ldots,a_{w+1,m}\right)
\ \ \ \ \ \ \ \ \ \ \left(  \text{since }\mathbf{d}_{u}\left(  w\right)
=w+1\right)  .
\end{align*}
Comparing this with%
\begin{align*}
A_{w+1,\bullet}  &  =\left(  \text{the }\left(  w+1\right)  \text{-th row of
the matrix }A\right) \\
&  \ \ \ \ \ \ \ \ \ \ \left(
\begin{array}
[c]{c}%
\text{since }A_{w+1,\bullet}\text{ is the }\left(  w+1\right)  \text{-th row
of the matrix }A\\
\text{(by the definition of }A_{w+1,\bullet}\text{)}%
\end{array}
\right) \\
&  =\left(  a_{w+1,1},a_{w+1,2},\ldots,a_{w+1,m}\right)
\ \ \ \ \ \ \ \ \ \ \left(  \text{since }A=\left(  a_{i,j}\right)  _{1\leq
i\leq n,\ 1\leq j\leq m}\right)  ,
\end{align*}
we obtain $\left(  A_{\sim u,\bullet}\right)  _{w,\bullet}=A_{w+1,\bullet}$.
This proves Proposition \ref{prop.unrows.basics} \textbf{(g)}.

\textbf{(h)} Let $v\in\left\{  1,2,\ldots,m\right\}  $ and $w\in\left\{
1,2,\ldots,v-1\right\}  $. From $w\in\left\{  1,2,\ldots,v-1\right\}  $, we
obtain $w\geq1$ and $w\leq v-1$. Thus, $w\leq v-1<v$, so that $w\leq v\leq m$
(since $v\in\left\{  1,2,\ldots,m\right\}  $). Also, $w\in\left\{
1,2,\ldots,v-1\right\}  \subseteq\left\{  1,2,\ldots,m-1\right\}  $ (since
$\underbrace{v}_{\leq m}-1\leq m-1$). From $w\geq1$ and $w\leq m$, we obtain
$w\in\left\{  1,2,\ldots,m\right\}  $. Thus, Proposition
\ref{prop.sol.unrows.2} \textbf{(b)} (applied to $m$, $w$ and $v$ instead of
$n$, $r$ and $s$) shows that
\begin{equation}
\mathbf{d}_{w,v}\left(  j\right)  =\mathbf{d}_{v}\left(  \mathbf{d}_{w}\left(
j\right)  \right)  \ \ \ \ \ \ \ \ \ \ \text{for every }j\in\mathbb{Z}.
\label{pf.prop.unrows.basics.h.1}%
\end{equation}

Now, $A=\left(  a_{i,j}\right)  _{1\leq i\leq n,\ 1\leq j\leq m}$. Hence,
Proposition \ref{prop.sol.unirows.d} \textbf{(b)} yields
\[
A_{\bullet,\sim v}=\left(  a_{i,\mathbf{d}_{v}\left(  j\right)  }\right)
_{1\leq i\leq n,\ 1\leq j\leq m-1}.
\]
Thus, Proposition \ref{prop.sol.unirows.d} \textbf{(b)} (applied to $m-1$,
$A_{\bullet,\sim v}$, $a_{i,\mathbf{d}_{v}\left(  j\right)  }$ and $w$ instead
of $m$, $A$, $a_{i,j}$ and $v$) yields
\begin{align}
\left(  A_{\bullet,\sim v}\right)  _{\bullet,\sim w}  &  =\left(
a_{i,\mathbf{d}_{v}\left(  \mathbf{d}_{w}\left(  j\right)  \right)  }\right)
_{1\leq i\leq n,\ 1\leq j\leq\left(  m-1\right)  -1}=\left(
\underbrace{a_{i,\mathbf{d}_{v}\left(  \mathbf{d}_{w}\left(  j\right)
\right)  }}_{\substack{=a_{i,\mathbf{d}_{w,v}\left(  j\right)  }\\\text{(since
}\mathbf{d}_{v}\left(  \mathbf{d}_{w}\left(  j\right)  \right)  =\mathbf{d}%
_{w,v}\left(  j\right)  \\\text{(by (\ref{pf.prop.unrows.basics.h.1})))}%
}}\right)  _{1\leq i\leq n,\ 1\leq j\leq m-2}\nonumber\\
&  \ \ \ \ \ \ \ \ \ \ \left(  \text{since }\left(  m-1\right)  -1=m-2\right)
\nonumber\\
&  =\left(  a_{i,\mathbf{d}_{w,v}\left(  j\right)  }\right)  _{1\leq i\leq
n,\ 1\leq j\leq m-2}. \label{pf.prop.unrows.basics.h.3}%
\end{align}
On the other hand, recall that $w$ and $v$ are two elements of $\left\{
1,2,\ldots,m\right\}  $ and satisfy $w<v$. Hence, Proposition
\ref{prop.sol.unrows.2} \textbf{(a)} (applied to $m$, $w$ and $v$ instead of
$n$, $r$ and $s$) yields%
\[
\left(  1,2,\ldots,\widehat{w},\ldots,\widehat{v},\ldots,m\right)  =\left(
\mathbf{d}_{w,v}\left(  1\right)  ,\mathbf{d}_{w,v}\left(  2\right)
,\ldots,\mathbf{d}_{w,v}\left(  m-2\right)  \right)  .
\]
Thus,%
\begin{align*}
\operatorname*{cols}\nolimits_{1,2,\ldots,\widehat{w},\ldots,\widehat{v}%
,\ldots,m}A  &  =\operatorname*{cols}\nolimits_{\mathbf{d}_{w,v}\left(
1\right)  ,\mathbf{d}_{w,v}\left(  2\right)  ,\ldots,\mathbf{d}_{w,v}\left(
m-2\right)  }A=\left(  a_{i,\mathbf{d}_{w,v}\left(  y\right)  }\right)
_{1\leq i\leq n,\ 1\leq y\leq m-2}\\
&  \ \ \ \ \ \ \ \ \ \ \left(
\begin{array}
[c]{c}%
\text{by the definition of }\operatorname*{cols}\nolimits_{\mathbf{d}%
_{w,v}\left(  1\right)  ,\mathbf{d}_{w,v}\left(  2\right)  ,\ldots
,\mathbf{d}_{w,v}\left(  m-2\right)  }A\\
\text{(since }A=\left(  a_{i,j}\right)  _{1\leq i\leq n,\ 1\leq j\leq
m}\text{)}%
\end{array}
\right) \\
&  =\left(  a_{i,\mathbf{d}_{w,v}\left(  j\right)  }\right)  _{1\leq i\leq
n,\ 1\leq j\leq m-2}%
\end{align*}
(here, we have renamed the index $\left(  i,y\right)  $ as $\left(
i,j\right)  $). Comparing this with (\ref{pf.prop.unrows.basics.h.3}), we
obtain $\left(  A_{\bullet,\sim v}\right)  _{\bullet,\sim w}%
=\operatorname*{cols}\nolimits_{1,2,\ldots,\widehat{w},\ldots,\widehat{v}%
,\ldots,m}A$. This proves Proposition \ref{prop.unrows.basics} \textbf{(h)}.

\textbf{(i)} Let $v\in\left\{  1,2,\ldots,m\right\}  $ and $w\in\left\{
v,v+1,\ldots,m-1\right\}  $. From $w\in\left\{  v,v+1,\ldots,m-1\right\}  $,
we obtain $w\geq v$ and $w\leq m-1$. Now, $w+1>w\geq v$, so that $w+1\geq
v\geq1$ (since $v\in\left\{  1,2,\ldots,m\right\}  $). Combining this with
$w+1\leq m$ (since $w\leq m-1$), we obtain $w+1\in\left\{  1,2,\ldots
,m\right\}  $. Also, $w\in\left\{  v,v+1,\ldots,m-1\right\}  \subseteq\left\{
1,2,\ldots,m-1\right\}  $ (since $v\geq1$). Now, $w+1>v$, so that $v<w+1$.
Thus, Proposition \ref{prop.sol.unrows.2} \textbf{(c)} (applied to $m$, $v$
and $w+1$ instead of $n$, $r$ and $s$) shows that%
\[
\mathbf{d}_{v,w+1}\left(  j\right)  =\mathbf{d}_{v}\left(  \mathbf{d}_{\left(
w+1\right)  -1}\left(  j\right)  \right)  \ \ \ \ \ \ \ \ \ \ \text{for every
}j\in\mathbb{Z}.
\]
In other words,
\begin{equation}
\mathbf{d}_{v,w+1}\left(  j\right)  =\mathbf{d}_{v}\left(  \mathbf{d}%
_{w}\left(  j\right)  \right)  \ \ \ \ \ \ \ \ \ \ \text{for every }%
j\in\mathbb{Z} \label{pf.prop.unrows.basics.i.1}%
\end{equation}
(since $\left(  w+1\right)  -1=w$).

Now, $A=\left(  a_{i,j}\right)  _{1\leq i\leq n,\ 1\leq j\leq m}$. Hence,
Proposition \ref{prop.sol.unirows.d} \textbf{(b)} yields
\[
A_{\bullet,\sim v}=\left(  a_{i,\mathbf{d}_{v}\left(  j\right)  }\right)
_{1\leq i\leq n,\ 1\leq j\leq m-1}.
\]
Thus, Proposition \ref{prop.sol.unirows.d} \textbf{(b)} (applied to $m-1$,
$A_{\bullet,\sim v}$, $a_{i,\mathbf{d}_{v}\left(  j\right)  }$ and $w$ instead
of $m$, $A$, $a_{i,j}$ and $v$) yields
\begin{align}
\left(  A_{\bullet,\sim v}\right)  _{\bullet,\sim w}  &  =\left(
a_{i,\mathbf{d}_{v}\left(  \mathbf{d}_{w}\left(  j\right)  \right)  }\right)
_{1\leq i\leq n,\ 1\leq j\leq\left(  m-1\right)  -1}=\left(
\underbrace{a_{i,\mathbf{d}_{v}\left(  \mathbf{d}_{w}\left(  j\right)
\right)  }}_{\substack{=a_{i,\mathbf{d}_{v,w+1}\left(  j\right)
}\\\text{(since }\mathbf{d}_{v}\left(  \mathbf{d}_{w}\left(  j\right)
\right)  =\mathbf{d}_{v,w+1}\left(  j\right)  \\\text{(by
(\ref{pf.prop.unrows.basics.i.1})))}}}\right)  _{1\leq i\leq n,\ 1\leq j\leq
m-2}\nonumber\\
&  \ \ \ \ \ \ \ \ \ \ \left(  \text{since }\left(  m-1\right)  -1=m-2\right)
\nonumber\\
&  =\left(  a_{i,\mathbf{d}_{v,w+1}\left(  j\right)  }\right)  _{1\leq i\leq
n,\ 1\leq j\leq m-2}. \label{pf.prop.unrows.basics.i.3}%
\end{align}

But $v$ and $w+1$ are two elements of $\left\{  1,2,\ldots,m\right\}  $ and
satisfy $v<w+1$. Hence, Proposition \ref{prop.sol.unrows.2} \textbf{(a)}
(applied to $m$, $v$ and $w+1$ instead of $n$, $r$ and $s$) yields%
\[
\left(  1,2,\ldots,\widehat{v},\ldots,\widehat{w+1},\ldots,m\right)  =\left(
\mathbf{d}_{v,w+1}\left(  1\right)  ,\mathbf{d}_{v,w+1}\left(  2\right)
,\ldots,\mathbf{d}_{v,w+1}\left(  m-2\right)  \right)  .
\]
Thus,%
\begin{align*}
\operatorname*{cols}\nolimits_{1,2,\ldots,\widehat{v},\ldots,\widehat{w+1}%
,\ldots,m}A  &  =\operatorname*{cols}\nolimits_{\mathbf{d}_{v,w+1}\left(
1\right)  ,\mathbf{d}_{v,w+1}\left(  2\right)  ,\ldots,\mathbf{d}%
_{v,w+1}\left(  m-2\right)  }A\\
&  =\left(  a_{i,\mathbf{d}_{v,w+1}\left(  y\right)  }\right)  _{1\leq i\leq
n,\ 1\leq y\leq m-2}\\
&  \ \ \ \ \ \ \ \ \ \ \left(
\begin{array}
[c]{c}%
\text{by the definition of }\operatorname*{cols}\nolimits_{\mathbf{d}%
_{v,w+1}\left(  1\right)  ,\mathbf{d}_{v,w+1}\left(  2\right)  ,\ldots
,\mathbf{d}_{v,w+1}\left(  m-2\right)  }A\\
\text{(since }A=\left(  a_{i,j}\right)  _{1\leq i\leq n,\ 1\leq j\leq
m}\text{)}%
\end{array}
\right) \\
&  =\left(  a_{i,\mathbf{d}_{v,w+1}\left(  j\right)  }\right)  _{1\leq i\leq
n,\ 1\leq j\leq m-2}%
\end{align*}
(here, we have renamed the index $\left(  i,y\right)  $ as $\left(
i,j\right)  $). Comparing this with (\ref{pf.prop.unrows.basics.i.3}), we
obtain $\left(  A_{\bullet,\sim v}\right)  _{\bullet,\sim w}%
=\operatorname*{cols}\nolimits_{1,2,\ldots,\widehat{v},\ldots,\widehat{w+1}%
,\ldots,m}A$. This proves Proposition \ref{prop.unrows.basics} \textbf{(i)}.

\textbf{(j)} Let $u\in\left\{  1,2,\ldots,n\right\}  $ and $w\in\left\{
1,2,\ldots,u-1\right\}  $. From $w\in\left\{  1,2,\ldots,u-1\right\}  $, we
obtain $w\geq1$ and $w\leq u-1$. Thus, $w\leq u-1<u$, so that $w\leq u\leq n$
(since $u\in\left\{  1,2,\ldots,n\right\}  $). Also, $w\in\left\{
1,2,\ldots,u-1\right\}  \subseteq\left\{  1,2,\ldots,n-1\right\}  $ (since
$\underbrace{u}_{\leq n}-1\leq n-1$). From $w\geq1$ and $w\leq n$, we obtain
$w\in\left\{  1,2,\ldots,n\right\}  $. Thus, Proposition
\ref{prop.sol.unrows.2} \textbf{(b)} (applied to $m$, $w$ and $u$ instead of
$n$, $r$ and $s$) shows that
\[
\mathbf{d}_{w,u}\left(  j\right)  =\mathbf{d}_{u}\left(  \mathbf{d}_{w}\left(
j\right)  \right)  \ \ \ \ \ \ \ \ \ \ \text{for every }j\in\mathbb{Z}.
\]
Renaming the index $j$ as $i$ in this result, we obtain the following fact:%
\begin{equation}
\mathbf{d}_{w,u}\left(  i\right)  =\mathbf{d}_{u}\left(  \mathbf{d}_{w}\left(
i\right)  \right)  \ \ \ \ \ \ \ \ \ \ \text{for every }i\in\mathbb{Z}.
\label{pf.prop.unrows.basics.j.1}%
\end{equation}

Now, $A=\left(  a_{i,j}\right)  _{1\leq i\leq n,\ 1\leq j\leq m}$. Hence,
Proposition \ref{prop.sol.unirows.d} \textbf{(a)} yields
\[
A_{\sim u,\bullet}=\left(  a_{\mathbf{d}_{u}\left(  i\right)  ,j}\right)
_{1\leq i\leq n-1,\ 1\leq j\leq m}.
\]
Thus, Proposition \ref{prop.sol.unirows.d} \textbf{(a)} (applied to $n-1$,
$A_{\sim u,\bullet}$, $a_{\mathbf{d}_{u}\left(  i\right)  ,j}$ and $w$ instead
of $n$, $A$, $a_{i,j}$ and $u$) yields
\begin{align}
\left(  A_{\sim u,\bullet}\right)  _{\sim w,\bullet}  &  =\left(
a_{\mathbf{d}_{u}\left(  \mathbf{d}_{w}\left(  i\right)  \right)  ,j}\right)
_{1\leq i\leq\left(  n-1\right)  -1,\ 1\leq j\leq m}=\left(
\underbrace{a_{\mathbf{d}_{u}\left(  \mathbf{d}_{w}\left(  i\right)  \right)
,j}}_{\substack{=a_{\mathbf{d}_{w,u}\left(  i\right)  ,j}\\\text{(since
}\mathbf{d}_{u}\left(  \mathbf{d}_{w}\left(  i\right)  \right)  =\mathbf{d}%
_{w,u}\left(  i\right)  \\\text{(by (\ref{pf.prop.unrows.basics.j.1})))}%
}}\right)  _{1\leq i\leq n-2,\ 1\leq j\leq m}\nonumber\\
&  \ \ \ \ \ \ \ \ \ \ \left(  \text{since }\left(  n-1\right)  -1=n-2\right)
\nonumber\\
&  =\left(  a_{\mathbf{d}_{w,u}\left(  i\right)  ,j}\right)  _{1\leq i\leq
n-2,\ 1\leq j\leq m}. \label{pf.prop.unrows.basics.j.3}%
\end{align}
On the other hand, recall that $w$ and $u$ are two elements of $\left\{
1,2,\ldots,n\right\}  $ and satisfy $w<u$. Hence, Proposition
\ref{prop.sol.unrows.2} \textbf{(a)} (applied to $w$ and $u$ instead of $r$
and $s$) yields%
\[
\left(  1,2,\ldots,\widehat{w},\ldots,\widehat{u},\ldots,n\right)  =\left(
\mathbf{d}_{w,u}\left(  1\right)  ,\mathbf{d}_{w,u}\left(  2\right)
,\ldots,\mathbf{d}_{w,u}\left(  n-2\right)  \right)  .
\]
Thus,%
\begin{align*}
\operatorname*{rows}\nolimits_{1,2,\ldots,\widehat{w},\ldots,\widehat{u}%
,\ldots,n}A  &  =\operatorname*{rows}\nolimits_{\mathbf{d}_{w,u}\left(
1\right)  ,\mathbf{d}_{w,u}\left(  2\right)  ,\ldots,\mathbf{d}_{w,u}\left(
n-2\right)  }A=\left(  a_{\mathbf{d}_{w,u}\left(  x\right)  ,j}\right)
_{1\leq x\leq n-2,\ 1\leq j\leq m}\\
&  \ \ \ \ \ \ \ \ \ \ \left(  \text{by the definition of }%
\operatorname*{rows}\nolimits_{\mathbf{d}_{w,v}\left(  1\right)
,\mathbf{d}_{w,v}\left(  2\right)  ,\ldots,\mathbf{d}_{w,v}\left(  n-2\right)
}A\right) \\
&  =\left(  a_{\mathbf{d}_{w,u}\left(  i\right)  ,j}\right)  _{1\leq i\leq
n-2,\ 1\leq j\leq m}%
\end{align*}
(here, we have renamed the index $\left(  x,j\right)  $ as $\left(
i,j\right)  $). Comparing this with (\ref{pf.prop.unrows.basics.j.3}), we
obtain $\left(  A_{\sim u,\bullet}\right)  _{\sim w,\bullet}%
=\operatorname*{rows}\nolimits_{1,2,\ldots,\widehat{w},\ldots,\widehat{u}%
,\ldots,n}A$. This proves Proposition \ref{prop.unrows.basics} \textbf{(j)}.

\textbf{(k)} Let $u\in\left\{  1,2,\ldots,n\right\}  $ and $w\in\left\{
u,u+1,\ldots,n-1\right\}  $. From $w\in\left\{  u,u+1,\ldots,n-1\right\}  $,
we obtain $w\geq u$ and $w\leq n-1$. But $u\geq1$ (since $u\in\left\{
1,2,\ldots,n\right\}  $). Hence, $w\geq u\geq1$. Therefore, $w+1\geq w\geq1$.
Combining this with $w+1\leq n$ (since $w\leq n-1$), we obtain $w+1\in\left\{
1,2,\ldots,n\right\}  $. Also, $w\in\left\{  u,u+1,\ldots,n-1\right\}
\subseteq\left\{  1,2,\ldots,n-1\right\}  $ (since $u\geq1$). Now, $w+1>w\geq
u$, so that $u<w+1$. Thus, Proposition \ref{prop.sol.unrows.2} \textbf{(c)}
(applied to $u$ and $w+1$ instead of $r$ and $s$) shows that%
\[
\mathbf{d}_{u,w+1}\left(  j\right)  =\mathbf{d}_{u}\left(  \mathbf{d}_{\left(
w+1\right)  -1}\left(  j\right)  \right)  \ \ \ \ \ \ \ \ \ \ \text{for every
}j\in\mathbb{Z}.
\]
In other words,
\[
\mathbf{d}_{u,w+1}\left(  j\right)  =\mathbf{d}_{u}\left(  \mathbf{d}%
_{w}\left(  j\right)  \right)  \ \ \ \ \ \ \ \ \ \ \text{for every }%
j\in\mathbb{Z}%
\]
(since $\left(  w+1\right)  -1=w$). Renaming the index $j$ as $i$ in this
result, we obtain the following fact:%
\begin{equation}
\mathbf{d}_{u,w+1}\left(  i\right)  =\mathbf{d}_{u}\left(  \mathbf{d}%
_{w}\left(  i\right)  \right)  \ \ \ \ \ \ \ \ \ \ \text{for every }%
i\in\mathbb{Z}. \label{pf.prop.unrows.basics.k.1}%
\end{equation}

Now, $A=\left(  a_{i,j}\right)  _{1\leq i\leq n,\ 1\leq j\leq m}$. Hence,
Proposition \ref{prop.sol.unirows.d} \textbf{(a)} yields
\[
A_{\sim u,\bullet}=\left(  a_{\mathbf{d}_{u}\left(  i\right)  ,j}\right)
_{1\leq i\leq n-1,\ 1\leq j\leq m}.
\]
Thus, Proposition \ref{prop.sol.unirows.d} \textbf{(a)} (applied to $n-1$,
$A_{\sim u,\bullet}$, $a_{\mathbf{d}_{u}\left(  i\right)  ,j}$ and $w$ instead
of $n$, $A$, $a_{i,j}$ and $u$) yields
\begin{align}
\left(  A_{\sim u,\bullet}\right)  _{\sim w,\bullet}  &  =\left(
a_{\mathbf{d}_{u}\left(  \mathbf{d}_{w}\left(  i\right)  \right)  ,j}\right)
_{1\leq i\leq\left(  n-1\right)  -1,\ 1\leq j\leq m}\nonumber\\
&  =\left(  \underbrace{a_{\mathbf{d}_{u}\left(  \mathbf{d}_{w}\left(
i\right)  \right)  ,j}}_{\substack{=a_{\mathbf{d}_{u,w+1}\left(  i\right)
,j}\\\text{(since }\mathbf{d}_{u}\left(  \mathbf{d}_{w}\left(  i\right)
\right)  =\mathbf{d}_{u,w+1}\left(  i\right)  \\\text{(by
(\ref{pf.prop.unrows.basics.k.1})))}}}\right)  _{1\leq i\leq n-2,\ 1\leq j\leq
m}\nonumber\\
&  \ \ \ \ \ \ \ \ \ \ \left(  \text{since }\left(  n-1\right)  -1=n-2\right)
\nonumber\\
&  =\left(  a_{\mathbf{d}_{u,w+1}\left(  i\right)  ,j}\right)  _{1\leq i\leq
n-2,\ 1\leq j\leq m}. \label{pf.prop.unrows.basics.k.3}%
\end{align}
On the other hand, recall that $u$ and $w+1$ are two elements of $\left\{
1,2,\ldots,n\right\}  $ and satisfy $u<w+1$. Hence, Proposition
\ref{prop.sol.unrows.2} \textbf{(a)} (applied to $u$ and $w+1$ instead of $r$
and $s$) yields%
\[
\left(  1,2,\ldots,\widehat{u},\ldots,\widehat{w+1},\ldots,n\right)  =\left(
\mathbf{d}_{u,w+1}\left(  1\right)  ,\mathbf{d}_{u,w+1}\left(  2\right)
,\ldots,\mathbf{d}_{u,w+1}\left(  n-2\right)  \right)  .
\]
Thus,%
\begin{align*}
\operatorname*{rows}\nolimits_{1,2,\ldots,\widehat{u},\ldots,\widehat{w+1}%
,\ldots,n}A  &  =\operatorname*{rows}\nolimits_{\mathbf{d}_{u,w+1}\left(
1\right)  ,\mathbf{d}_{u,w+1}\left(  2\right)  ,\ldots,\mathbf{d}%
_{u,w+1}\left(  n-2\right)  }A\\
&  =\left(  a_{\mathbf{d}_{u,w+1}\left(  x\right)  ,j}\right)  _{1\leq x\leq
n-2,\ 1\leq j\leq m}\\
&  \ \ \ \ \ \ \ \ \ \ \left(  \text{by the definition of }%
\operatorname*{rows}\nolimits_{\mathbf{d}_{u,w+1}\left(  1\right)
,\mathbf{d}_{u,w+1}\left(  2\right)  ,\ldots,\mathbf{d}_{u,w+1}\left(
n-2\right)  }A\right) \\
&  =\left(  a_{\mathbf{d}_{u,w+1}\left(  i\right)  ,j}\right)  _{1\leq i\leq
n-2,\ 1\leq j\leq m}%
\end{align*}
(here, we have renamed the index $\left(  x,j\right)  $ as $\left(
i,j\right)  $). Comparing this with (\ref{pf.prop.unrows.basics.k.3}), we
obtain $\left(  A_{\sim u,\bullet}\right)  _{\sim w,\bullet}%
=\operatorname*{rows}\nolimits_{1,2,\ldots,\widehat{u},\ldots,\widehat{w+1}%
,\ldots,n}A$. This proves Proposition \ref{prop.unrows.basics} \textbf{(k)}.

\textbf{(l)} Let $v\in\left\{  1,2,\ldots,n\right\}  $, $u\in\left\{
1,2,\ldots,n\right\}  $ and $q\in\left\{  1,2,\ldots,m\right\}  $ be such that
$u<v$. We have $u<v$ and thus $u\leq v-1$ (since $u$ and $v$ are integers).
Combined with $u\geq1$ (since $u\in\left\{  1,2,\ldots,n\right\}  $), this
shows that $u\in\left\{  1,2,\ldots,v-1\right\}  $.

Proposition \ref{prop.unrows.basics} \textbf{(c)} (applied to $v$ and $q$
instead of $u$ and $v$) yields $\left(  A_{\bullet,\sim q}\right)  _{\sim
v,\bullet}=\left(  A_{\sim v,\bullet}\right)  _{\bullet,\sim q}=A_{\sim v,\sim
q}$.

We have $v\in\left\{  1,2,\ldots,n\right\}  $, thus $v\leq n$. Now,
$u\leq\underbrace{v}_{\leq n}-1\leq n-1$. Combining this with $u\geq1$, we
obtain $u\in\left\{  1,2,\ldots,n-1\right\}  $. Hence, Proposition
\ref{prop.unrows.basics} \textbf{(c)} (applied to $n-1$, $A_{\sim v,\bullet}$
and $q$ instead of $n$, $A$ and $v$) yields $\left(  \left(  A_{\sim
v,\bullet}\right)  _{\bullet,\sim q}\right)  _{\sim u,\bullet}=\left(  \left(
A_{\sim v,\bullet}\right)  _{\sim u,\bullet}\right)  _{\bullet,\sim q}=\left(
A_{\sim v,\bullet}\right)  _{\sim u,\sim q}$. Hence,%
\begin{align*}
\left(  A_{\sim v,\bullet}\right)  _{\sim u,\sim q}  &  =\left(
\underbrace{\left(  A_{\sim v,\bullet}\right)  _{\bullet,\sim q}}_{=\left(
A_{\bullet,\sim q}\right)  _{\sim v,\bullet}}\right)  _{\sim u,\bullet
}=\left(  \left(  A_{\bullet,\sim q}\right)  _{\sim v,\bullet}\right)  _{\sim
u,\bullet}\\
&  =\operatorname*{rows}\nolimits_{1,2,\ldots,\widehat{u},\ldots
,\widehat{v},\ldots,n}\left(  A_{\bullet,\sim q}\right)
\end{align*}
(by Proposition \ref{prop.unrows.basics} \textbf{(j)}, applied to $m-1$,
$A_{\bullet,\sim q}$, $v$ and $u$ instead of $m$, $A$, $u$ and $w$). This
proves Proposition \ref{prop.unrows.basics} \textbf{(l)}.
\end{proof}

\begin{proof}
[Proof of Proposition \ref{prop.unrows.basics-I}.]We have $u<v\leq n$ (since
$v\in\left\{  1,2,\ldots,n\right\}  $) and thus $u\leq n-1$ (since $u$ and $n$
are integers). Combining this with $u\geq1$ (since $u\in\left\{
1,2,\ldots,n\right\}  $), we obtain $u\in\left\{  1,2,\ldots,n-1\right\}  $.
Hence, the column vector $\left(  I_{n-1}\right)  _{\bullet,u}\in
\mathbb{K}^{\left(  n-1\right)  \times1}$ is well-defined. Also, the column
vector $\left(  I_{n}\right)  _{\bullet,u}\in\mathbb{K}^{n\times1}$ is
well-defined (since $u\in\left\{  1,2,\ldots,n\right\}  $), and therefore the
column vector $\left(  \left(  I_{n}\right)  _{\bullet,u}\right)  _{\sim
v,\bullet}\in\mathbb{K}^{\left(  n-1\right)  \times1}$ is well-defined.

For any two objects $i$ and $j$, we define an element $\delta_{i,j}%
\in\mathbb{K}$ by $\delta_{i,j}=%
\begin{cases}
1, & \text{if }i=j;\\
0, & \text{if }i\neq j
\end{cases}
$. Then, every $m\in\mathbb{N}$ satisfies%
\begin{equation}
I_{m}=\left(  \delta_{i,j}\right)  _{1\leq i\leq m,\ 1\leq j\leq m}
\label{pf.prop.unrows.basics-I.Im}%
\end{equation}
(by the definition of $I_{m}$). Applying this to $m=n$, we obtain
$I_{n}=\left(  \delta_{i,j}\right)  _{1\leq i\leq n,\ 1\leq j\leq n}$. Now,
$\left(  I_{n}\right)  _{\bullet,u}$ is the $u$-th column of the matrix
$I_{n}$ (by the definition of $\left(  I_{n}\right)  _{\bullet,u}$). Thus,%
\begin{align*}
\left(  I_{n}\right)  _{\bullet,u}  &  =\left(  \text{the }u\text{-th column
of the matrix }I_{n}\right) \\
&  =\left(
\begin{array}
[c]{c}%
\delta_{1,u}\\
\delta_{2,u}\\
\vdots\\
\delta_{n,u}%
\end{array}
\right)  \ \ \ \ \ \ \ \ \ \ \left(  \text{since }I_{n}=\left(  \delta
_{i,j}\right)  _{1\leq i\leq n,\ 1\leq j\leq n}\right) \\
&  =\left(  \delta_{i,u}\right)  _{1\leq i\leq n,\ 1\leq j\leq1}.
\end{align*}
Hence, Proposition \ref{prop.sol.unirows.d} \textbf{(a)} (applied to $1$,
$\left(  I_{n}\right)  _{\bullet,u}$, $\delta_{i,u}$ and $v$ instead of $m$,
$A$, $a_{i,j}$ and $u$) yields%
\begin{equation}
\left(  \left(  I_{n}\right)  _{\bullet,u}\right)  _{\sim v,\bullet}=\left(
\delta_{\mathbf{d}_{v}\left(  i\right)  ,u}\right)  _{1\leq i\leq n-1,\ 1\leq
j\leq1}. \label{pf.prop.unrows.basics-I.1}%
\end{equation}

But every $i\in\mathbb{Z}$ satisfies%
\begin{equation}
\delta_{\mathbf{d}_{v}\left(  i\right)  ,u}=\delta_{i,u}
\label{pf.prop.unrows.basics-I.2}%
\end{equation}
\footnote{\textit{Proof of (\ref{pf.prop.unrows.basics-I.1}):} Let
$i\in\mathbb{Z}$. We must prove the equality (\ref{pf.prop.unrows.basics-I.2}%
).
\par
The definition of $\mathbf{d}_{v}\left(  i\right)  $ yields $\mathbf{d}%
_{v}\left(  i\right)  =%
\begin{cases}
i, & \text{if }i<v;\\
i+1, & \text{if }i\geq v
\end{cases}
$. But $u<v$, and thus $v>u$.
\par
We are in one of the following two Cases:
\par
\textit{Case 1:} We have $i<v$.
\par
\textit{Case 2:} We have $i\geq v$.
\par
Let us first consider Case 1. In this case, we have $i<v$. Now, $\mathbf{d}%
_{v}\left(  i\right)  =%
\begin{cases}
i, & \text{if }i<v;\\
i+1, & \text{if }i\geq v
\end{cases}
=i$ (since $i<v$) and thus $\delta_{\mathbf{d}_{v}\left(  i\right)  ,u}%
=\delta_{i,u}$. Hence, (\ref{pf.prop.unrows.basics-I.2}) is proven in Case 1.
\par
Let us now consider Case 2. In this case, we have $i\geq v$. Now,
$\mathbf{d}_{v}\left(  i\right)  =%
\begin{cases}
i, & \text{if }i<v;\\
i+1, & \text{if }i\geq v
\end{cases}
=i+1$ (since $i\geq v$) and thus $\mathbf{d}_{v}\left(  i\right)  =i+1>i\geq
v>u$. Hence, $\mathbf{d}_{v}\left(  i\right)  \neq u$. Thus, $\delta
_{\mathbf{d}_{v}\left(  i\right)  ,u}=0$. Also, $i\geq v>u$ and thus $i\neq
u$; hence, $\delta_{i,u}=0$. Now, $\delta_{\mathbf{d}_{v}\left(  i\right)
,u}=0=\delta_{i,u}$. Hence, (\ref{pf.prop.unrows.basics-I.2}) is proven in
Case 2.
\par
We have now proven (\ref{pf.prop.unrows.basics-I.2}) in each of the two Cases
1 and 2. Since these two Cases cover all possibilities, this shows that
(\ref{pf.prop.unrows.basics-I.2}) always holds.}. Hence,
(\ref{pf.prop.unrows.basics-I.1}) becomes%
\begin{equation}
\left(  \left(  I_{n}\right)  _{\bullet,u}\right)  _{\sim v,\bullet}=\left(
\underbrace{\delta_{\mathbf{d}_{v}\left(  i\right)  ,u}}_{\substack{=\delta
_{i,u}\\\text{(by (\ref{pf.prop.unrows.basics-I.2}))}}}\right)  _{1\leq i\leq
n-1,\ 1\leq j\leq1}=\left(  \delta_{i,u}\right)  _{1\leq i\leq n-1,\ 1\leq
j\leq1}. \label{pf.prop.unrows.basics-I.3}%
\end{equation}

But (\ref{pf.prop.unrows.basics-I.Im}) (applied to $m=n-1$) yields
$I_{n-1}=\left(  \delta_{i,j}\right)  _{1\leq i\leq n-1,\ 1\leq j\leq n-1}$.
Now, $\left(  I_{n-1}\right)  _{\bullet,u}$ is the $u$-th column of the matrix
$I_{n-1}$ (by the definition of $\left(  I_{n-1}\right)  _{\bullet,u}$). Thus,%
\begin{align*}
\left(  I_{n-1}\right)  _{\bullet,u}  &  =\left(  \text{the }u\text{-th column
of the matrix }I_{n-1}\right) \\
&  =\left(
\begin{array}
[c]{c}%
\delta_{1,u}\\
\delta_{2,u}\\
\vdots\\
\delta_{n-1,u}%
\end{array}
\right)  \ \ \ \ \ \ \ \ \ \ \left(  \text{since }I_{n-1}=\left(  \delta
_{i,j}\right)  _{1\leq i\leq n-1,\ 1\leq j\leq n-1}\right) \\
&  =\left(  \delta_{i,u}\right)  _{1\leq i\leq n-1,\ 1\leq j\leq1}.
\end{align*}
Comparing this with (\ref{pf.prop.unrows.basics-I.3}), we obtain $\left(
\left(  I_{n}\right)  _{\bullet,u}\right)  _{\sim v,\bullet}=\left(
I_{n-1}\right)  _{\bullet,u}$. This proves Proposition
\ref{prop.unrows.basics-I}.
\end{proof}

We have now proven Proposition \ref{prop.unrows.basics} and Proposition
\ref{prop.unrows.basics-I}. Thus, Exercise \ref{exe.unrows.basics}
\textbf{(a)} is solved.

Next, let us derive Proposition \ref{prop.desnanot.12} and Proposition
\ref{prop.desnanot.1n} from Theorem \ref{thm.desnanot}:

\begin{proof}
[Proof of Proposition \ref{prop.desnanot.12} using Theorem \ref{thm.desnanot}%
.]We have%
\begin{equation}
\operatorname*{sub}\nolimits_{1,2,\ldots,\widehat{1},\ldots,\widehat{2}%
,\ldots,n}^{1,2,\ldots,\widehat{1},\ldots,\widehat{2},\ldots,n}A=\widetilde{A}
\label{pf.prop.desnanot.12.1}%
\end{equation}
\footnote{\textit{Proof of (\ref{pf.prop.desnanot.12.1}):} The definition of
$\left(  1,2,\ldots,\widehat{1},\ldots,\widehat{2},\ldots,n\right)  $ yields%
\begin{align*}
&  \left(  1,2,\ldots,\widehat{1},\ldots,\widehat{2},\ldots,n\right) \\
&  =\left(  \underbrace{1,2,\ldots,1-1}_{\substack{\text{all integers}%
\\\text{from }1\text{ to }1-1}},\underbrace{1+1,1+2,\ldots,2-1}%
_{\substack{\text{all integers}\\\text{from }1+1\text{ to }2-1}%
},\underbrace{2+1,2+2,\ldots,n}_{\substack{\text{all integers}\\\text{from
}2+1\text{ to }n}}\right) \\
&  =\left(  \underbrace{\underbrace{1,2,\ldots,0}_{\substack{\text{all
integers}\\\text{from }1\text{ to }0}}}_{\text{this is an empty list}%
},\underbrace{\underbrace{2,3,\ldots,1}_{\substack{\text{all integers}%
\\\text{from }2\text{ to }1}}}_{\text{this is an empty list}}%
,\underbrace{3,4,\ldots,n}_{\substack{\text{all integers}\\\text{from }3\text{
to }n}}\right) \\
&  =\left(  \underbrace{3,4,\ldots,n}_{\substack{\text{all integers}%
\\\text{from }3\text{ to }n}}\right)  =\left(  3,4,\ldots,n\right)  =\left(
1+2,2+2,\ldots,\left(  n-2\right)  +2\right)  .
\end{align*}
Hence,%
\begin{align*}
\operatorname*{sub}\nolimits_{1,2,\ldots,\widehat{1},\ldots,\widehat{2}%
,\ldots,n}^{1,2,\ldots,\widehat{1},\ldots,\widehat{2},\ldots,n}A  &
=\operatorname*{sub}\nolimits_{1+2,2+2,\ldots,\left(  n-2\right)
+2}^{1+2,2+2,\ldots,\left(  n-2\right)  +2}A=\left(  a_{x+2,y+2}\right)
_{1\leq x\leq n-2,\ 1\leq y\leq n-2}\\
&  \ \ \ \ \ \ \ \ \ \ \left(
\begin{array}
[c]{c}%
\text{by the definition of }\operatorname*{sub}\nolimits_{1+2,2+2,\ldots
,\left(  n-2\right)  +2}^{1+2,2+2,\ldots,\left(  n-2\right)  +2}A\text{,}\\
\text{since }A=\left(  a_{i,j}\right)  _{1\leq i\leq n,\ 1\leq j\leq n}%
\end{array}
\right) \\
&  =\left(  a_{i+2,j+2}\right)  _{1\leq i\leq n-2,\ 1\leq j\leq n-2}\\
&  \ \ \ \ \ \ \ \ \ \ \left(  \text{here, we have renamed the index }\left(
x,y\right)  \text{ as }\left(  i,j\right)  \right) \\
&  =\widetilde{A}\ \ \ \ \ \ \ \ \ \ \left(  \text{since }\widetilde{A}%
=\left(  a_{i+2,j+2}\right)  _{1\leq i\leq n-2,\ 1\leq j\leq n-2}\right)  .
\end{align*}
Qed.}. Now, $1<2$ and $1<2$. Moreover, $1$ and $2$ are elements of $\left\{
1,2,\ldots,n\right\}  $ (since $n\geq2$). Hence, Theorem \ref{thm.desnanot}
(applied to $p=1$, $q=2$, $u=1$ and $v=2$) yields%
\begin{align*}
&  \det A\cdot\det\left(  \operatorname*{sub}\nolimits_{1,2,\ldots
,\widehat{1},\ldots,\widehat{2},\ldots,n}^{1,2,\ldots,\widehat{1}%
,\ldots,\widehat{2},\ldots,n}A\right) \\
&  =\det\left(  A_{\sim1,\sim1}\right)  \cdot\det\left(  A_{\sim2,\sim
2}\right)  -\det\left(  A_{\sim1,\sim2}\right)  \cdot\det\left(  A_{\sim
2,\sim1}\right)  .
\end{align*}
In view of (\ref{pf.prop.desnanot.12.1}), this rewrites as
\begin{align*}
&  \det A\cdot\det\widetilde{A}\\
&  =\det\left(  A_{\sim1,\sim1}\right)  \cdot\det\left(  A_{\sim2,\sim
2}\right)  -\det\left(  A_{\sim1,\sim2}\right)  \cdot\det\left(  A_{\sim
2,\sim1}\right)  .
\end{align*}
This proves Proposition \ref{prop.desnanot.12}.
\end{proof}

\begin{proof}
[Proof of Proposition \ref{prop.desnanot.1n} using Theorem \ref{thm.desnanot}%
.]We have%
\begin{equation}
\operatorname*{sub}\nolimits_{1,2,\ldots,\widehat{1},\ldots,\widehat{n}%
,\ldots,n}^{1,2,\ldots,\widehat{1},\ldots,\widehat{n},\ldots,n}A=A^{\prime}
\label{pf.prop.desnanot.1n.1}%
\end{equation}
\footnote{\textit{Proof of (\ref{pf.prop.desnanot.1n.1}):} The definition of
$\left(  1,2,\ldots,\widehat{1},\ldots,\widehat{n},\ldots,n\right)  $ yields%
\begin{align*}
&  \left(  1,2,\ldots,\widehat{1},\ldots,\widehat{n},\ldots,n\right) \\
&  =\left(  \underbrace{\underbrace{1,2,\ldots,1-1}_{\substack{\text{all
integers}\\\text{from }1\text{ to }1-1}}}_{\text{this is an empty list}%
},\underbrace{1+1,1+2,\ldots,n-1}_{\substack{\text{all integers}\\\text{from
}1+1\text{ to }n-1}},\underbrace{\underbrace{n+1,n+2,\ldots,n}%
_{\substack{\text{all integers}\\\text{from }n+1\text{ to }n}}}_{\text{this is
an empty list}}\right) \\
&  =\left(  \underbrace{1+1,1+2,\ldots,n-1}_{\substack{\text{all
integers}\\\text{from }1+1\text{ to }n-1}}\right)  =\left(  1+1,1+2,\ldots
,n-1\right) \\
&  =\left(  2,3,\ldots,n-1\right)  =\left(  1+1,2+1,\ldots,\left(  n-2\right)
+1\right)  .
\end{align*}
Hence,%
\begin{align*}
\operatorname*{sub}\nolimits_{1,2,\ldots,\widehat{1},\ldots,\widehat{n}%
,\ldots,n}^{1,2,\ldots,\widehat{1},\ldots,\widehat{n},\ldots,n}A  &
=\operatorname*{sub}\nolimits_{1+1,2+1,\ldots,\left(  n-2\right)
+1}^{1+1,2+1,\ldots,\left(  n-2\right)  +1}A=\left(  a_{x+1,y+1}\right)
_{1\leq x\leq n-2,\ 1\leq y\leq n-2}\\
&  \ \ \ \ \ \ \ \ \ \ \left(
\begin{array}
[c]{c}%
\text{by the definition of }\operatorname*{sub}\nolimits_{1+1,2+1,\ldots
,\left(  n-2\right)  +1}^{1+1,2+1,\ldots,\left(  n-2\right)  +1}A\text{,}\\
\text{since }A=\left(  a_{i,j}\right)  _{1\leq i\leq n,\ 1\leq j\leq n}%
\end{array}
\right) \\
&  =\left(  a_{i+1,j+1}\right)  _{1\leq i\leq n-2,\ 1\leq j\leq n-2}\\
&  \ \ \ \ \ \ \ \ \ \ \left(  \text{here, we have renamed the index }\left(
x,y\right)  \text{ as }\left(  i,j\right)  \right) \\
&  =A^{\prime}\ \ \ \ \ \ \ \ \ \ \left(  \text{since }A^{\prime}=\left(
a_{i+1,j+1}\right)  _{1\leq i\leq n-2,\ 1\leq j\leq n-2}\right)  .
\end{align*}
Qed.}. Now, $1<n$ (since $n\geq2$) and $1<n$. Moreover, $1$ and $n$ are
elements of $\left\{  1,2,\ldots,n\right\}  $ (since $n\geq2\geq1$). Hence,
Theorem \ref{thm.desnanot} (applied to $p=1$, $q=n$, $u=1$ and $v=n$) yields%
\begin{align*}
&  \det A\cdot\det\left(  \operatorname*{sub}\nolimits_{1,2,\ldots
,\widehat{1},\ldots,\widehat{n},\ldots,n}^{1,2,\ldots,\widehat{1}%
,\ldots,\widehat{n},\ldots,n}A\right) \\
&  =\det\left(  A_{\sim1,\sim1}\right)  \cdot\det\left(  A_{\sim n,\sim
n}\right)  -\det\left(  A_{\sim1,\sim n}\right)  \cdot\det\left(  A_{\sim
n,\sim1}\right)  .
\end{align*}
In view of (\ref{pf.prop.desnanot.1n.1}), this rewrites as
\begin{align*}
&  \det A\cdot\det\left(  A^{\prime}\right) \\
&  =\det\left(  A_{\sim1,\sim1}\right)  \cdot\det\left(  A_{\sim n,\sim
n}\right)  -\det\left(  A_{\sim1,\sim n}\right)  \cdot\det\left(  A_{\sim
n,\sim1}\right)  .
\end{align*}
This proves Proposition \ref{prop.desnanot.1n}.
\end{proof}

\begin{proof}
[Solution to Exercise \ref{exe.unrows.basics}.]\textbf{(a)} We proven
Proposition \ref{prop.unrows.basics} and Proposition
\ref{prop.unrows.basics-I} above. Thus, Exercise \ref{exe.unrows.basics}
\textbf{(a)} is solved.

\textbf{(b)} We have derived Proposition \ref{prop.desnanot.12} and
Proposition \ref{prop.desnanot.1n} from Theorem \ref{thm.desnanot} above.
Thus, Exercise \ref{exe.unrows.basics} \textbf{(b)} is solved.
\end{proof}
\end{verlong}

\subsection{Solution to Exercise \ref{exe.prop.addcol.props}}

\begin{vershort}
\begin{proof}
[Proof of Proposition \ref{prop.addcol.props1}.]The matrix $\left(  A\mid
v\right)  $ is defined as the $n\times\left(  m+1\right)  $-matrix whose $m+1$
columns are $A_{\bullet,1},A_{\bullet,2},\ldots,A_{\bullet,m},v$ (from left to
right). Thus, the first $m$ columns of the matrix $\left(  A\mid v\right)  $
are $A_{\bullet,1},A_{\bullet,2},\ldots,A_{\bullet,m}$, whereas the $\left(
m+1\right)  $-st column of $\left(  A\mid v\right)  $ is $v$.

\textbf{(a)} For every $q\in\left\{  1,2,\ldots,m\right\}  $, we have%
\begin{align*}
\left(  A\mid v\right)  _{\bullet,q}  &  =\left(  \text{the }q\text{-th column
of the matrix }\left(  A\mid v\right)  \right) \\
&  \ \ \ \ \ \ \ \ \ \ \left(  \text{by the definition of }\left(  A\mid
v\right)  _{\bullet,q}\right) \\
&  =A_{\bullet,q}%
\end{align*}
(since the first $m$ columns of the matrix $\left(  A\mid v\right)  $ are
$A_{\bullet,1},A_{\bullet,2},\ldots,A_{\bullet,m}$). This proves Proposition
\ref{prop.addcol.props1} \textbf{(a)}.

\textbf{(b)} The definition of $\left(  A\mid v\right)  _{\bullet,m+1}$ yields%
\[
\left(  A\mid v\right)  _{\bullet,m+1}=\left(  \text{the }\left(  m+1\right)
\text{-st column of the matrix }\left(  A\mid v\right)  \right)  =v
\]
(since the $\left(  m+1\right)  $-st column of $\left(  A\mid v\right)  $ is
$v$). This proves Proposition \ref{prop.addcol.props1} \textbf{(b)}.

\textbf{(c)} Let $q\in\left\{  1,2,\ldots,m\right\}  $. Recall that $\left(
A\mid v\right)  $ is the $n\times\left(  m+1\right)  $-matrix whose $m+1$
columns are $A_{\bullet,1},A_{\bullet,2},\ldots,A_{\bullet,m},v$. The matrix
$\left(  A\mid v\right)  _{\bullet,\sim q}$ is obtained from this matrix
$\left(  A\mid v\right)  $ by removing its $q$-th column; thus,
\begin{equation}
\text{the columns of this matrix }\left(  A\mid v\right)  _{\bullet,\sim
q}\text{ are }A_{\bullet,1},A_{\bullet,2},\ldots,\widehat{A_{\bullet,q}%
},\ldots,A_{\bullet,m},v \label{pf.prop.addcol.props1.short.c.1}%
\end{equation}
(since the columns of the matrix $\left(  A\mid v\right)  $ are $A_{\bullet
,1},A_{\bullet,2},\ldots,A_{\bullet,m},v$, and since $q\in\left\{
1,2,\ldots,m\right\}  $).

On the other hand, the matrix $A_{\bullet,\sim q}$ is obtained from the matrix
$A$ by removing its $q$-th column; thus, the columns of this matrix
$A_{\bullet,\sim q}$ are $A_{\bullet,1},A_{\bullet,2},\ldots
,\widehat{A_{\bullet,q}},\ldots,A_{\bullet,m}$. Hence,
\begin{equation}
\text{the columns of the matrix }\left(  A_{\bullet,\sim q}\mid v\right)
\text{ are }A_{\bullet,1},A_{\bullet,2},\ldots,\widehat{A_{\bullet,q}}%
,\ldots,A_{\bullet,m},v \label{pf.prop.addcol.props1.short.c.2}%
\end{equation}
(since the matrix $\left(  A_{\bullet,\sim q}\mid v\right)  $ is obtained from
$A_{\bullet,\sim q}$ by attaching the column $v$ at its \textquotedblleft
right edge\textquotedblright).

Comparing (\ref{pf.prop.addcol.props1.short.c.1}) with
(\ref{pf.prop.addcol.props1.short.c.2}), we see that the columns of the matrix
$\left(  A\mid v\right)  _{\bullet,\sim q}$ are precisely the columns of the
matrix $\left(  A_{\bullet,\sim q}\mid v\right)  $. Thus, these two matrices
must be identical. In other words, $\left(  A\mid v\right)  _{\bullet,\sim
q}=\left(  A_{\bullet,\sim q}\mid v\right)  $. This proves Proposition
\ref{prop.addcol.props1} \textbf{(c)}.

\textbf{(d)} We have defined $\left(  A\mid v\right)  $ as the $n\times\left(
m+1\right)  $-matrix whose $m+1$ columns are $A_{\bullet,1},A_{\bullet
,2},\ldots,A_{\bullet,m},v$. The matrix $\left(  A\mid v\right)
_{\bullet,\sim\left(  m+1\right)  }$ is obtained from this matrix $\left(
A\mid v\right)  $ by removing its $\left(  m+1\right)  $-th column; thus,
\begin{equation}
\text{the columns of this matrix }\left(  A\mid v\right)  _{\bullet
,\sim\left(  m+1\right)  }\text{ are }A_{\bullet,1},A_{\bullet,2}%
,\ldots,A_{\bullet,m} \label{pf.prop.addcol.props1.short.d.1}%
\end{equation}
(since the columns of the matrix $\left(  A\mid v\right)  $ are $A_{\bullet
,1},A_{\bullet,2},\ldots,A_{\bullet,m},v$). On the other hand, clearly,%
\begin{equation}
\text{the columns of the matrix }A\text{ are }A_{\bullet,1},A_{\bullet
,2},\ldots,A_{\bullet,m}. \label{pf.prop.addcol.props1.short.d.2}%
\end{equation}
Comparing (\ref{pf.prop.addcol.props1.short.d.1}) with
(\ref{pf.prop.addcol.props1.short.d.2}), we see that the columns of the matrix
$\left(  A\mid v\right)  _{\bullet,\sim\left(  m+1\right)  }$ are precisely
the columns of the matrix $A$. Thus, these two matrices must be identical. In
other words, $\left(  A\mid v\right)  _{\bullet,\sim\left(  m+1\right)  }=A$.
This proves Proposition \ref{prop.addcol.props1} \textbf{(d)}.

\textbf{(e)} Here is an informal proof: Let $p\in\left\{  1,2,\ldots
,n\right\}  $. The matrix $\left(  A\mid v\right)  $ is obtained from the
matrix $A$ by attaching the column $v$ at its \textquotedblleft right
edge\textquotedblright. The matrix $\left(  A\mid v\right)  _{\sim p,\bullet}$
is obtained from the matrix $\left(  A\mid v\right)  $ by removing its $p$-th
row. Hence, the matrix $\left(  A\mid v\right)  _{\sim p,\bullet}$ is obtained
from the matrix $A$ by first attaching the column $v$ at its \textquotedblleft
right edge\textquotedblright\ and then removing the $p$-th row. On the other
hand, the matrix $\left(  A_{\sim p,\bullet}\mid v_{\sim p,\bullet}\right)  $
is obtained from the matrix $A$ by first removing the $p$-th row (this is how
we get $A_{\sim p,\bullet}$), and then attaching the column $v_{\sim
p,\bullet}$ (this is $v$ with its $p$-th row removed) at its \textquotedblleft
right edge\textquotedblright. Obviously, the two procedures result in the same
matrix, i.e., we have $\left(  A\mid v\right)  _{\sim p,\bullet}=\left(
A_{\sim p,\bullet}\mid v_{\sim p,\bullet}\right)  $. Thus, Proposition
\ref{prop.addcol.props1} \textbf{(e)} is proven.\footnote{Here is an outline
of a more formal proof: First, we observe that if $A$ is an $n\times
m$-matrix, if $p\in\left\{  1,2,\ldots,n\right\}  $ and if $q\in\left\{
1,2,\ldots,m\right\}  $, then%
\begin{equation}
\left(  A_{\sim p,\bullet}\right)  _{\bullet,q}=\left(  A_{\bullet,q}\right)
_{\sim p,\bullet}. \label{pf.prop.addcol.props1.short.e.fn1.pf.1}%
\end{equation}
(Indeed, this is a simple fact that is similar to the statements of
Proposition \ref{prop.unrows.basics}. What it says is that if you remove the
$p$-th row of the matrix $A$ and then take the $q$-th column of the resulting
matrix, then you get the same result as when you first take the $q$-th column
of $A$ and then remove the $p$-th row of this column.)
\par
Now, let us prove Proposition \ref{prop.addcol.props1} \textbf{(e)} formally:
Both $\left(  A\mid v\right)  _{\sim p,\bullet}$ and $\left(  A_{\sim
p,\bullet}\mid v_{\sim p,\bullet}\right)  $ are $\left(  n-1\right)
\times\left(  m+1\right)  $-matrices (since they have one fewer row and one
more column than $A$). We shall now show that
\begin{equation}
\left(  \left(  A\mid v\right)  _{\sim p,\bullet}\right)  _{\bullet,q}=\left(
A_{\sim p,\bullet}\mid v_{\sim p,\bullet}\right)  _{\bullet,q}
\label{pf.prop.addcol.props1.short.e.fn1.pf.8}%
\end{equation}
for every $q\in\left\{  1,2,\ldots,m+1\right\}  $.
\par
\textit{Proof of (\ref{pf.prop.addcol.props1.short.e.fn1.pf.8}):} Let
$q\in\left\{  1,2,\ldots,m+1\right\}  $. We must prove
(\ref{pf.prop.addcol.props1.short.e.fn1.pf.8}).
\par
We can apply (\ref{pf.prop.addcol.props1.short.e.fn1.pf.1}) to $m+1$ and
$\left(  A\mid v\right)  $ instead of $m$ and $A$. As a result, we obtain%
\begin{equation}
\left(  \left(  A\mid v\right)  _{\sim p,\bullet}\right)  _{\bullet,q}=\left(
\left(  A\mid v\right)  _{\bullet,q}\right)  _{\sim p,\bullet}.
\label{pf.prop.addcol.props1.short.e.fn1.pf.8.pf.1}%
\end{equation}
\par
We are in one of the following two cases:
\par
\textit{Case 1:} We have $q\neq m+1$.
\par
\textit{Case 2:} We have $q=m+1$.
\par
Let us consider Case 1 first. In this case, we have $q\neq m+1$. Hence,
$q\in\left\{  1,2,\ldots,m\right\}  $ (since $q\in\left\{  1,2,\ldots
,m+1\right\}  $). Thus, Proposition \ref{prop.addcol.props1} \textbf{(a)}
yields $\left(  A\mid v\right)  _{\bullet,q}=A_{\bullet,q}$. Now,
(\ref{pf.prop.addcol.props1.short.e.fn1.pf.8.pf.1}) becomes%
\[
\left(  \left(  A\mid v\right)  _{\sim p,\bullet}\right)  _{\bullet,q}=\left(
\underbrace{\left(  A\mid v\right)  _{\bullet,q}}_{=A_{\bullet,q}}\right)
_{\sim p,\bullet}=\left(  A_{\bullet,q}\right)  _{\sim p,\bullet}.
\]
Comparing this with%
\begin{align*}
\left(  A_{\sim p,\bullet}\mid v_{\sim p,\bullet}\right)  _{\bullet,q}  &
=\left(  A_{\sim p,\bullet}\right)  _{\bullet,q}\ \ \ \ \ \ \ \ \ \ \left(
\begin{array}
[c]{c}%
\text{by Proposition \ref{prop.addcol.props1} \textbf{(a)}, applied to}\\
n-1\text{, }A_{\sim p,\bullet}\text{ and }v_{\sim p,\bullet}\text{ instead of
}n\text{, }A\text{ and }v
\end{array}
\right) \\
&  =\left(  A_{\bullet,q}\right)  _{\sim p,\bullet}\ \ \ \ \ \ \ \ \ \ \left(
\text{by (\ref{pf.prop.addcol.props1.short.e.fn1.pf.1})}\right)  ,
\end{align*}
we obtain $\left(  \left(  A\mid v\right)  _{\sim p,\bullet}\right)
_{\bullet,q}=\left(  A_{\sim p,\bullet}\mid v_{\sim p,\bullet}\right)
_{\bullet,q}$. Thus, (\ref{pf.prop.addcol.props1.short.e.fn1.pf.8}) is proven
in Case 1.
\par
Let us now consider Case 2. In this case, we have $q=m+1$. Hence, $\left(
A\mid v\right)  _{\bullet,q}=\left(  A\mid v\right)  _{\bullet,m+1}=v$ (by
Proposition \ref{prop.addcol.props1} \textbf{(b)}). Now,
(\ref{pf.prop.addcol.props1.short.e.fn1.pf.8.pf.1}) becomes%
\[
\left(  \left(  A\mid v\right)  _{\sim p,\bullet}\right)  _{\bullet,q}=\left(
\underbrace{\left(  A\mid v\right)  _{\bullet,q}}_{=v}\right)  _{\sim
p,\bullet}=v_{\sim p,\bullet}.
\]
Comparing this with%
\begin{align*}
\left(  A_{\sim p,\bullet}\mid v_{\sim p,\bullet}\right)  _{\bullet,q}  &
=\left(  A_{\sim p,\bullet}\mid v_{\sim p,\bullet}\right)  _{\bullet
,m+1}\ \ \ \ \ \ \ \ \ \ \left(  \text{since }q=m+1\right) \\
&  =v_{\sim p,\bullet}\ \ \ \ \ \ \ \ \ \ \left(
\begin{array}
[c]{c}%
\text{by Proposition \ref{prop.addcol.props1} \textbf{(b)}, applied to}\\
n-1\text{, }A_{\sim p,\bullet}\text{ and }v_{\sim p,\bullet}\text{ instead of
}n\text{, }A\text{ and }v
\end{array}
\right)  ,
\end{align*}
we obtain $\left(  \left(  A\mid v\right)  _{\sim p,\bullet}\right)
_{\bullet,q}=\left(  A_{\sim p,\bullet}\mid v_{\sim p,\bullet}\right)
_{\bullet,q}$. Thus, (\ref{pf.prop.addcol.props1.short.e.fn1.pf.8}) is proven
in Case 2.
\par
Now we have proven (\ref{pf.prop.addcol.props1.short.e.fn1.pf.8}) in each of
the two Cases 1 and 2. Thus, the proof of
(\ref{pf.prop.addcol.props1.short.e.fn1.pf.8}) is complete.
\par
Now, for every $q\in\left\{  1,2,\ldots,m+1\right\}  $, we have%
\begin{align*}
&  \left(  \text{the }q\text{-th column of the matrix }\left(  A\mid v\right)
_{\sim p,\bullet}\right) \\
&  =\left(  \left(  A\mid v\right)  _{\sim p,\bullet}\right)  _{\bullet,q}\\
&  =\left(  A_{\sim p,\bullet}\mid v_{\sim p,\bullet}\right)  _{\bullet
,q}\ \ \ \ \ \ \ \ \ \ \left(  \text{by
(\ref{pf.prop.addcol.props1.short.e.fn1.pf.8})}\right) \\
&  =\left(  \text{the }q\text{-th column of the matrix }\left(  A_{\sim
p,\bullet}\mid v_{\sim p,\bullet}\right)  \right)  .
\end{align*}
This shows that the matrices $\left(  A\mid v\right)  _{\sim p,\bullet}$ and
$\left(  A_{\sim p,\bullet}\mid v_{\sim p,\bullet}\right)  $ are identical
(since they are both $\left(  n-1\right)  \times\left(  m+1\right)
$-matrices). Proposition \ref{prop.addcol.props1} \textbf{(e)} is thus
proven.}

\textbf{(f)} Let $p\in\left\{  1,2,\ldots,n\right\}  $. Then, $\left(  A\mid
v\right)  $ is an $n\times\left(  m+1\right)  $-matrix. Proposition
\ref{prop.unrows.basics} \textbf{(c)} (applied to $m+1$, $\left(  A\mid
v\right)  $, $p$ and $m+1$ instead of $m$, $A$, $u$ and $v$) yields
\[
\left(  \left(  A\mid v\right)  _{\bullet,\sim\left(  m+1\right)  }\right)
_{\sim p,\bullet}=\left(  \left(  A\mid v\right)  _{\sim p,\bullet}\right)
_{\bullet,\sim\left(  m+1\right)  }=\left(  A\mid v\right)  _{\sim
p,\sim\left(  m+1\right)  }.
\]
Hence,%
\[
\left(  A\mid v\right)  _{\sim p,\sim\left(  m+1\right)  }=\left(
\underbrace{\left(  A\mid v\right)  _{\bullet,\sim\left(  m+1\right)  }%
}_{\substack{=A\\\text{(by Proposition \ref{prop.addcol.props1} \textbf{(d)}%
)}}}\right)  _{\sim p,\bullet}=A_{\sim p,\bullet}.
\]
This proves Proposition \ref{prop.addcol.props1} \textbf{(f)}.
\end{proof}

\begin{proof}
[Proof of Proposition \ref{prop.addcol.props2}.]\textbf{(a)} Let $v=\left(
v_{1},v_{2},\ldots,v_{n}\right)  ^{T}\in\mathbb{K}^{n\times1}$ be a vector.

Clearly, $\left(  A\mid v\right)  $ is an $n\times n$-matrix (since it is
obtained by attaching the column $v$ to the $n\times\left(  n-1\right)
$-matrix $A$). Write it in the form $\left(  A\mid v\right)  =\left(
b_{i,j}\right)  _{1\leq i\leq n,\ 1\leq j\leq n}$. Then, every $p\in\left\{
1,2,\ldots,n\right\}  $ satisfies%
\begin{equation}
b_{p,n}=v_{p} \label{pf.prop.addcols.props2.short.1}%
\end{equation}
\footnote{\textit{Proof of (\ref{pf.prop.addcols.props2.short.1}):} We have
$\left(  A\mid v\right)  =\left(  b_{i,j}\right)  _{1\leq i\leq n,\ 1\leq
j\leq n}$. Thus,%
\[
\left(  \text{the }n\text{-th column of the matrix }\left(  A\mid v\right)
\right)  =\left(
\begin{array}
[c]{c}%
b_{1,n}\\
b_{2,n}\\
\vdots\\
b_{n,n}%
\end{array}
\right)  .
\]
Hence,%
\begin{align*}
\left(
\begin{array}
[c]{c}%
b_{1,n}\\
b_{2,n}\\
\vdots\\
b_{n,n}%
\end{array}
\right)   &  =\left(  \text{the }n\text{-th column of the matrix }\left(
A\mid v\right)  \right)  =\left(  A\mid v\right)  _{\bullet,n}\\
&  =\left(  A\mid v\right)  _{\bullet,\left(  n-1\right)  +1}%
\ \ \ \ \ \ \ \ \ \ \left(  \text{since }n=\left(  n-1\right)  +1\right) \\
&  =v\ \ \ \ \ \ \ \ \ \ \left(  \text{by Proposition \ref{prop.addcol.props1}
\textbf{(b)}, applied to }m=n-1\right) \\
&  =\left(  v_{1},v_{2},\ldots,v_{n}\right)  ^{T}=\left(
\begin{array}
[c]{c}%
v_{1}\\
v_{2}\\
\vdots\\
v_{n}%
\end{array}
\right)  .
\end{align*}
In other words, every $p\in\left\{  1,2,\ldots,n\right\}  $ satisfies
$b_{p,n}=v_{p}$. This proves (\ref{pf.prop.addcols.props2.short.1}).} and%
\begin{equation}
\left(  A\mid v\right)  _{\sim p,\sim n}=A_{\sim p,\bullet}
\label{pf.prop.addcols.props2.short.2}%
\end{equation}
\footnote{\textit{Proof of (\ref{pf.prop.addcols.props2.short.2}):} Let
$p\in\left\{  1,2,\ldots,n\right\}  $. Proposition \ref{prop.addcol.props1}
\textbf{(f)} (applied to $m=n-1$) yields $\left(  A\mid v\right)  _{\sim
p,\sim\left(  \left(  n-1\right)  +1\right)  }=A_{\sim p,\bullet}$. This
rewrites as $\left(  A\mid v\right)  _{\sim p,\sim n}=A_{\sim p,\bullet}$
(since $\left(  n-1\right)  +1=n$). This proves
(\ref{pf.prop.addcols.props2.short.2}).}.

We have $n\in\left\{  1,2,\ldots,n\right\}  $ (since $n$ is a positive
integer) and $\left(  A\mid v\right)  =\left(  b_{i,j}\right)  _{1\leq i\leq
n,\ 1\leq j\leq n}$. Hence, Theorem \ref{thm.laplace.gen} \textbf{(b)}
(applied to $\left(  A\mid v\right)  $, $b_{i,j}$ and $n$ instead of $A$,
$a_{i,j}$ and $q$) yields%
\begin{align*}
\det\left(  A\mid v\right)   &  =\sum_{p=1}^{n}\underbrace{\left(  -1\right)
^{p+n}}_{=\left(  -1\right)  ^{n+p}}\underbrace{b_{p,n}}_{\substack{=v_{p}%
\\\text{(by (\ref{pf.prop.addcols.props2.short.1}))}}}\det\left(
\underbrace{\left(  A\mid v\right)  _{\sim p,\sim n}}_{\substack{=A_{\sim
p,\bullet}\\\text{(by (\ref{pf.prop.addcols.props2.short.2}))}}}\right) \\
&  =\sum_{p=1}^{n}\left(  -1\right)  ^{n+p}v_{p}\det\left(  A_{\sim p,\bullet
}\right)  =\sum_{i=1}^{n}\left(  -1\right)  ^{n+i}v_{i}\det\left(  A_{\sim
i,\bullet}\right)
\end{align*}
(here, we have renamed the summation index $p$ as $i$). This proves
Proposition \ref{prop.addcol.props2} \textbf{(a)}.

\textbf{(b)} Let $p\in\left\{  1,2,\ldots,n\right\}  $. For any two objects
$i$ and $j$, we define an element $\delta_{i,j}\in\mathbb{K}$ by $\delta
_{i,j}=%
\begin{cases}
1, & \text{if }i=j;\\
0, & \text{if }i\neq j
\end{cases}
$. Then, $I_{n}=\left(  \delta_{i,j}\right)  _{1\leq i\leq n,\ 1\leq j\leq n}$
(by the definition of $I_{n}$). Now, $\left(  I_{n}\right)  _{\bullet,p}$ is
the $p$-th column of the matrix $I_{n}$ (by the definition of $\left(
I_{n}\right)  _{\bullet,p}$). Thus,%
\begin{align*}
\left(  I_{n}\right)  _{\bullet,p}  &  =\left(  \text{the }p\text{-th column
of the matrix }I_{n}\right) \\
&  =\left(
\begin{array}
[c]{c}%
\delta_{1,p}\\
\delta_{2,p}\\
\vdots\\
\delta_{n,p}%
\end{array}
\right)  \ \ \ \ \ \ \ \ \ \ \left(  \text{since }I_{n}=\left(  \delta
_{i,j}\right)  _{1\leq i\leq n,\ 1\leq j\leq n}\right) \\
&  =\left(  \delta_{1,p},\delta_{2,p},\ldots,\delta_{n,p}\right)  ^{T}.
\end{align*}
Hence, Proposition \ref{prop.addcol.props2} \textbf{(a)} (applied to $\left(
I_{n}\right)  _{\bullet,p}$ and $\delta_{i,p}$ instead of $v$ and $v_{i}$)
yields%
\begin{align*}
&  \det\left(  A\mid\left(  I_{n}\right)  _{\bullet,p}\right) \\
&  =\underbrace{\sum_{i=1}^{n}}_{=\sum_{i\in\left\{  1,2,\ldots,n\right\}  }%
}\left(  -1\right)  ^{n+i}\delta_{i,p}\det\left(  A_{\sim i,\bullet}\right)
=\sum_{i\in\left\{  1,2,\ldots,n\right\}  }\left(  -1\right)  ^{n+i}%
\delta_{i,p}\det\left(  A_{\sim i,\bullet}\right) \\
&  =\left(  -1\right)  ^{n+p}\underbrace{\delta_{p,p}}%
_{\substack{=1\\\text{(since }p=p\text{)}}}\det\left(  A_{\sim p,\bullet
}\right)  +\sum_{\substack{i\in\left\{  1,2,\ldots,n\right\}  ;\\i\neq
p}}\left(  -1\right)  ^{n+i}\underbrace{\delta_{i,p}}%
_{\substack{=0\\\text{(since }i\neq p\text{)}}}\det\left(  A_{\sim i,\bullet
}\right) \\
&  \ \ \ \ \ \ \ \ \ \ \left(
\begin{array}
[c]{c}%
\text{here, we have split off the addend for }i=p\text{ from the sum,}\\
\text{since }p\in\left\{  1,2,\ldots,n\right\}
\end{array}
\right) \\
&  =\left(  -1\right)  ^{n+p}\det\left(  A_{\sim p,\bullet}\right)
+\underbrace{\sum_{\substack{i\in\left\{  1,2,\ldots,n\right\}  ;\\i\neq
p}}\left(  -1\right)  ^{n+i}0\det\left(  A_{\sim i,\bullet}\right)  }%
_{=0}=\left(  -1\right)  ^{n+p}\det\left(  A_{\sim p,\bullet}\right)  .
\end{align*}
This proves Proposition \ref{prop.addcol.props2} \textbf{(b)}.
\end{proof}

Before we prove Proposition \ref{prop.addcol.props3}, let us make a simple observation:

\begin{lemma}
\label{lem.sol.prop.addcol.props3.bc.short}Let $n\in\mathbb{N}$. Let
$A=\left(  a_{i,j}\right)  _{1\leq i\leq n,\ 1\leq j\leq n}\in\mathbb{K}%
^{n\times n}$ be an $n\times n$-matrix. Let $r$ and $q$ be two elements of
$\left\{  1,2,\ldots,n\right\}  $. Then,%
\begin{equation}
\det\left(  A_{\bullet,\sim q}\mid A_{\bullet,r}\right)  =\left(  -1\right)
^{n+q}\sum_{p=1}^{n}\left(  -1\right)  ^{p+q}a_{p,r}\det\left(  A_{\sim p,\sim
q}\right)  .\nonumber
\end{equation}

\end{lemma}

\begin{proof}
[Proof of Lemma \ref{lem.sol.prop.addcol.props3.bc.short}.]We have
$r\in\left\{  1,2,\ldots,n\right\}  $. Hence, $A_{\bullet,r}$ is an $n\times
1$-matrix (since $A$ is an $n\times n$-matrix). Also, $r\in\left\{
1,2,\ldots,n\right\}  $ shows that $1\leq r\leq n$; hence, $1\leq n$, so that
$n$ is a positive integer.

Also, $q\in\left\{  1,2,\ldots,n\right\}  $. Hence, $A_{\bullet,\sim q}$ is an
$n\times\left(  n-1\right)  $-matrix (since $A$ is an $n\times n$-matrix).

Furthermore, the definition of $A_{\bullet,r}$ yields%
\begin{align*}
A_{\bullet,r}  &  =\left(  \text{the }r\text{-th column of the matrix
}A\right) \\
&  =\left(
\begin{array}
[c]{c}%
a_{1,r}\\
a_{2,r}\\
\vdots\\
a_{n,r}%
\end{array}
\right)  \ \ \ \ \ \ \ \ \ \ \left(  \text{since }A=\left(  a_{i,j}\right)
_{1\leq i\leq n,\ 1\leq j\leq n}\right) \\
&  =\left(  a_{1,r},a_{2,r},\ldots,a_{n,r}\right)  ^{T}.
\end{align*}
Thus, Proposition \ref{prop.addcol.props2} \textbf{(a)} (applied to
$A_{\bullet,\sim q}$, $A_{\bullet,r}$ and $a_{i,r}$ instead of $A$, $v$ and
$v_{i}$) yields%
\begin{equation}
\det\left(  A_{\bullet,\sim q}\mid A_{\bullet,r}\right)  =\sum_{i=1}%
^{n}\left(  -1\right)  ^{n+i}a_{i,r}\det\left(  \left(  A_{\bullet,\sim
q}\right)  _{\sim i,\bullet}\right)  .
\label{pf.lem.sol.prop.addcol.props3.bc.short.1}%
\end{equation}
But every $i\in\left\{  1,2,\ldots,n\right\}  $ satisfies%
\begin{equation}
\left(  A_{\bullet,\sim q}\right)  _{\sim i,\bullet}=A_{\sim i,\sim q}
\label{pf.lem.sol.prop.addcol.props3.bc.short.2}%
\end{equation}
\footnote{\textit{Proof of (\ref{pf.lem.sol.prop.addcol.props3.bc.short.2}):}
Let $i\in\left\{  1,2,\ldots,n\right\}  $. Then, Proposition
\ref{prop.unrows.basics} \textbf{(c)} (applied to $n$, $i$ and $q$ instead of
$m$, $u$ and $v$) yields $\left(  A_{\bullet,\sim q}\right)  _{\sim i,\bullet
}=\left(  A_{\sim i,\bullet}\right)  _{\bullet,\sim q}=A_{\sim i,\sim q}$.
This proves (\ref{pf.lem.sol.prop.addcol.props3.bc.short.2}).} and%
\begin{equation}
\left(  -1\right)  ^{n+i}=\left(  -1\right)  ^{n+q}\left(  -1\right)  ^{i+q}
\label{pf.lem.sol.prop.addcol.props3.bc.short.3}%
\end{equation}
\footnote{\textit{Proof of (\ref{pf.lem.sol.prop.addcol.props3.bc.short.3}):}
Let $i\in\left\{  1,2,\ldots,n\right\}  $. Then, $\left(  -1\right)
^{n+q}\left(  -1\right)  ^{i+q}=\left(  -1\right)  ^{\left(  n+q\right)
+\left(  i+q\right)  }=\left(  -1\right)  ^{n+i}$ (since $\left(  n+q\right)
+\left(  i+q\right)  =n+i+2q\equiv n+i\operatorname{mod}2$). This proves
(\ref{pf.lem.sol.prop.addcol.props3.bc.short.3}).}. Now,
(\ref{pf.lem.sol.prop.addcol.props3.bc.short.1}) becomes%
\begin{align*}
\det\left(  A_{\bullet,\sim q}\mid A_{\bullet,r}\right)   &  =\sum_{i=1}%
^{n}\underbrace{\left(  -1\right)  ^{n+i}}_{\substack{=\left(  -1\right)
^{n+q}\left(  -1\right)  ^{i+q}\\\text{(by
(\ref{pf.lem.sol.prop.addcol.props3.bc.short.3}))}}}a_{i,r}\det\left(
\underbrace{\left(  A_{\bullet,\sim q}\right)  _{\sim i,\bullet}%
}_{\substack{=A_{\sim i,\sim q}\\\text{(by
(\ref{pf.lem.sol.prop.addcol.props3.bc.short.2}))}}}\right) \\
&  =\sum_{i=1}^{n}\left(  -1\right)  ^{n+q}\left(  -1\right)  ^{i+q}%
a_{i,r}\det\left(  A_{\sim i,\sim q}\right) \\
&  =\left(  -1\right)  ^{n+q}\sum_{i=1}^{n}\left(  -1\right)  ^{i+q}%
a_{i,r}\det\left(  A_{\sim i,\sim q}\right) \\
&  =\left(  -1\right)  ^{n+q}\sum_{p=1}^{n}\left(  -1\right)  ^{p+q}%
a_{p,r}\det\left(  A_{\sim p,\sim q}\right) \\
&  \ \ \ \ \ \ \ \ \ \ \left(  \text{here, we have renamed the summation index
}i\text{ as }p\right)  .
\end{align*}
This proves Lemma \ref{lem.sol.prop.addcol.props3.bc.short}.
\end{proof}

\begin{proof}
[Proof of Proposition \ref{prop.addcol.props3}.]\textbf{(a)} Proposition
\ref{prop.addcol.props3} \textbf{(a)} simply says that if we remove the $n$-th
column from the matrix $A$, and then reattach this column back to the matrix
(at its right edge), then we get the original matrix $A$ back. This should be
completely obvious\footnote{Nevertheless, let me give a formal proof of
Proposition \ref{prop.addcol.props3} \textbf{(a)} as well:
\par
Assume that $n>0$. Thus, $n\in\left\{  1,2,\ldots,n\right\}  $. Hence,
$A_{\bullet,\sim n}$ is a well-defined $n\times\left(  n-1\right)  $-matrix,
and $A_{\bullet,n}$ is a well-defined $n\times1$-matrix. Therefore, $\left(
A_{\bullet,\sim n}\mid A_{\bullet,n}\right)  $ is an $n\times n$-matrix.
\par
We shall now show that%
\begin{equation}
A_{\bullet,q}=\left(  A_{\bullet,\sim n}\mid A_{\bullet,n}\right)
_{\bullet,q} \label{pf.prop.addcol.props3.short.a.main}%
\end{equation}
for each $q\in\left\{  1,2,\ldots,n\right\}  $.
\par
\textit{Proof of (\ref{pf.prop.addcol.props3.short.a.main}):} Let
$q\in\left\{  1,2,\ldots,n\right\}  $. We must prove the equality
(\ref{pf.prop.addcol.props3.short.a.main}). We are in one of the following two
cases:
\par
\textit{Case 1:} We have $q\neq n$.
\par
\textit{Case 2:} We have $q=n$.
\par
Let us first consider Case 1. In this case, we have $q\neq n$. Combining
$q\in\left\{  1,2,\ldots,n\right\}  $ with $q\neq n$, we obtain $q\in\left\{
1,2,\ldots,n\right\}  \setminus\left\{  n\right\}  =\left\{  1,2,\ldots
,n-1\right\}  $. Thus, Proposition \ref{prop.addcol.props1} \textbf{(a)}
(applied to $n-1$, $A_{\bullet,\sim n}$ and $A_{\bullet,n}$ instead of $m$,
$A$ and $v$) yields $\left(  A_{\bullet,\sim n}\mid A_{\bullet,n}\right)
_{\bullet,q}=\left(  A_{\bullet,\sim n}\right)  _{\bullet,q}$.
\par
But Proposition \ref{prop.unrows.basics} \textbf{(d)} (applied to $n$, $n$ and
$q$ instead of $m$, $v$ and $w$) yields $\left(  A_{\bullet,\sim n}\right)
_{\bullet,q}=A_{\bullet,q}$. Hence, $A_{\bullet,q}=\left(  A_{\bullet,\sim
n}\right)  _{\bullet,q}=\left(  A_{\bullet,\sim n}\mid A_{\bullet,n}\right)
_{\bullet,q}$ (since $\left(  A_{\bullet,\sim n}\mid A_{\bullet,n}\right)
_{\bullet,q}=\left(  A_{\bullet,\sim n}\right)  _{\bullet,q}$). Hence,
(\ref{pf.prop.addcol.props3.short.a.main}) is proven in Case 1.
\par
Let us now consider Case 2. In this case, we have $q=n$. But Proposition
\ref{prop.addcol.props1} \textbf{(b)} (applied to $n-1$, $A_{\bullet,\sim n}$
and $A_{\bullet,n}$ instead of $m$, $A$ and $v$) yields $\left(
A_{\bullet,\sim n}\mid A_{\bullet,n}\right)  _{\bullet,\left(  n-1\right)
+1}=A_{\bullet,n}$. Thus, $A_{\bullet,n}=\left(  A_{\bullet,\sim n}\mid
A_{\bullet,n}\right)  _{\bullet,\left(  n-1\right)  +1}=\left(  A_{\bullet
,\sim n}\mid A_{\bullet,n}\right)  _{\bullet,q}$ (since $\left(  n-1\right)
+1=n=q$). Now, $q=n$, so that $A_{\bullet,q}=A_{\bullet,n}=\left(
A_{\bullet,\sim n}\mid A_{\bullet,n}\right)  _{\bullet,q}$. Hence,
(\ref{pf.prop.addcol.props3.short.a.main}) is proven in Case 2.
\par
Now, (\ref{pf.prop.addcol.props3.short.a.main}) is proven in each of the two
Cases 1 and 2. This completes the proof of
(\ref{pf.prop.addcol.props3.short.a.main}).
\par
Now we know that $A$ and $\left(  A_{\bullet,\sim n}\mid A_{\bullet,n}\right)
$ are two $n\times n$-matrices, and that every $q\in\left\{  1,2,\ldots
,n\right\}  $ satisfies%
\begin{align*}
\left(  \text{the }q\text{-th column of the matrix }A\right)   &
=A_{\bullet,q}=\left(  A_{\bullet,\sim n}\mid A_{\bullet,n}\right)
_{\bullet,q}\ \ \ \ \ \ \ \ \ \ \left(  \text{by
(\ref{pf.prop.addcol.props3.short.a.main})}\right) \\
&  =\left(  \text{the }q\text{-th column of the matrix }\left(  A_{\bullet
,\sim n}\mid A_{\bullet,n}\right)  \right)  .
\end{align*}
Hence, the matrices $A$ and $\left(  A_{\bullet,\sim n}\mid A_{\bullet
,n}\right)  $ are identical. This proves Proposition \ref{prop.addcol.props3}
\textbf{(a)}.}.

\textbf{(b)} Write the matrix $A$ in the form $A=\left(  a_{i,j}\right)
_{1\leq i\leq n,\ 1\leq j\leq n}$. (This is possible since $A$ is an $n\times
n$-matrix.) Lemma \ref{lem.sol.prop.addcol.props3.bc.short} (applied to $r=q$)
yields%
\begin{equation}
\det\left(  A_{\bullet,\sim q}\mid A_{\bullet,q}\right)  =\left(  -1\right)
^{n+q}\underbrace{\sum_{p=1}^{n}\left(  -1\right)  ^{p+q}a_{p,q}\det\left(
A_{\sim p,\sim q}\right)  }_{\substack{=\det A\\\text{(by Theorem
\ref{thm.laplace.gen} \textbf{(b)})}}}=\left(  -1\right)  ^{n+q}\det
A.\nonumber
\end{equation}
This proves Proposition \ref{prop.addcol.props3} \textbf{(b)}.

[\textit{Remark:} Let me outline an alternative proof of Proposition
\ref{prop.addcol.props3} \textbf{(b)}: Let $q\in\left\{  1,2,\ldots,n\right\}
$. The matrix $\left(  A_{\bullet,\sim q}\mid A_{\bullet,q}\right)  $ is
obtained from the matrix $A$ by removing the $q$-th column and then
reattaching this column to the right end of the matrix. This procedure can be
replaced by the following procedure, which clearly leads to the same result:

\begin{itemize}
\item Swap the $q$-th column of $A$ with the $\left(  q+1\right)  $-th column;

\item then swap the $\left(  q+1\right)  $-th column of the resulting matrix
with the $\left(  q+2\right)  $-th column;

\item then swap the $\left(  q+2\right)  $-th column of the resulting matrix
with the $\left(  q+3\right)  $-th column;

\item and so on, finally swapping the $\left(  n-1\right)  $-st column of the
matrix with the $n$-th column.
\end{itemize}

But this latter procedure is a sequence of $n-q$ swaps of two columns. Each
such swap multiplies the determinant of the matrix by $-1$ (according to
Exercise \ref{exe.ps4.6} \textbf{(b)}). Thus, the whole procedure multiplies
the determinant of the matrix by $\left(  -1\right)  ^{n-q}=\left(  -1\right)
^{n+q}$ (since $n-q\equiv n+q\operatorname{mod}2$). Since this procedure takes
the matrix $A$ to the matrix $\left(  A_{\bullet,\sim q}\mid A_{\bullet
,q}\right)  $, we thus conclude that $\det\left(  A_{\bullet,\sim q}\mid
A_{\bullet,q}\right)  =\left(  -1\right)  ^{n+q}\det A$. This proves
Proposition \ref{prop.addcol.props3} \textbf{(b)} again.]

\textbf{(c)} Let $r$ and $q$ be two elements of $\left\{  1,2,\ldots
,n\right\}  $ satisfying $r\neq q$. We must show that $\det\left(
A_{\bullet,\sim q}\mid A_{\bullet,r}\right)  =0$.

Write the matrix $A$ in the form $A=\left(  a_{i,j}\right)  _{1\leq i\leq
n,\ 1\leq j\leq n}$. (This is possible since $A$ is an $n\times n$-matrix.)
Lemma \ref{lem.sol.prop.addcol.props3.bc.short} yields%
\begin{equation}
\det\left(  A_{\bullet,\sim q}\mid A_{\bullet,r}\right)  =\left(  -1\right)
^{n+q}\underbrace{\sum_{p=1}^{n}\left(  -1\right)  ^{p+q}a_{p,r}\det\left(
A_{\sim p,\sim q}\right)  }_{\substack{=0\\\text{(by Proposition
\ref{prop.laplace.0} \textbf{(b)}}\\\text{(since }q\neq r\text{))}%
}}=0.\nonumber
\end{equation}
This proves Proposition \ref{prop.addcol.props3} \textbf{(c)}.

[\textit{Remark:} Let me outline an alternative proof of Proposition
\ref{prop.addcol.props3} \textbf{(c)}: Let $r$ and $q$ be two elements of
$\left\{  1,2,\ldots,n\right\}  $ satisfying $r\neq q$. Then, $A_{\bullet,r}$
is one of the columns of the matrix $A_{\bullet,\sim q}$ (since the $r$-th
column of the matrix $A$ is not lost when the $q$-th column is removed).
Hence, the column vector $A_{\bullet,r}$ appears twice as a column in the
matrix $\left(  A_{\bullet,\sim q}\mid A_{\bullet,r}\right)  $ (namely, once
as one of the columns of $A_{\bullet,\sim q}$, and another time as the
attached column). Therefore, the matrix $\left(  A_{\bullet,\sim q}\mid
A_{\bullet,r}\right)  $ has two equal columns. Exercise \ref{exe.ps4.6}
\textbf{(f)} thus shows that $\det\left(  A_{\bullet,\sim q}\mid A_{\bullet
,r}\right)  =0$. This proves Proposition \ref{prop.addcol.props3} \textbf{(c)} again.]

\textbf{(d)} Let $p\in\left\{  1,2,\ldots,n\right\}  $ and $q\in\left\{
1,2,\ldots,n\right\}  $. Then, Proposition \ref{prop.unrows.basics}
\textbf{(c)} (applied to $u=p$ and $v=q$) yields%
\[
\left(  A_{\bullet,\sim q}\right)  _{\sim p,\bullet}=\left(  A_{\sim
p,\bullet}\right)  _{\bullet,\sim q}=A_{\sim p,\sim q}.
\]

But $A$ is an $n\times n$-matrix. Hence, $A_{\bullet,\sim q}$ is an
$n\times\left(  n-1\right)  $-matrix. Moreover, $p\in\left\{  1,2,\ldots
,n\right\}  $, so that $1\leq p\leq n$ and thus $n\geq1$; hence, $n$ is a
positive integer. Proposition \ref{prop.addcol.props2} \textbf{(b)} (applied
to $A_{\bullet,\sim q}$ instead of $A$) thus yields
\[
\det\left(  A_{\bullet,\sim q}\mid\left(  I_{n}\right)  _{\bullet,p}\right)
=\left(  -1\right)  ^{n+p}\det\left(  \underbrace{\left(  A_{\bullet,\sim
q}\right)  _{\sim p,\bullet}}_{=A_{\sim p,\sim q}}\right)  =\left(  -1\right)
^{n+p}\det\left(  A_{\sim p,\sim q}\right)  .
\]
This proves Proposition \ref{prop.addcol.props3} \textbf{(d)}.

\textbf{(e)} Let $u$ and $v$ be two elements of $\left\{  1,2,\ldots
,n\right\}  $ satisfying $u<v$. Let $r$ be an element of $\left\{
1,2,\ldots,n-1\right\}  $ satisfying $r\neq u$. We must show that $\det\left(
A_{\bullet,\sim u}\mid\left(  A_{\bullet,\sim v}\right)  _{\bullet,r}\right)
=0$.

We are in one of the following two cases:

\textit{Case 1:} We have $r<v$.

\textit{Case 2:} We have $r\geq v$.

Let us first consider Case 1. In this case, we have $r<v$. Thus, $r\in\left\{
1,2,\ldots,v-1\right\}  $. Hence, Proposition \ref{prop.unrows.basics}
\textbf{(d)} (applied to $m=n$ and $w=r$) yields $\left(  A_{\bullet,\sim
v}\right)  _{\bullet,r}=A_{\bullet,r}$. Hence,%
\[
\det\left(  A_{\bullet,\sim u}\mid\underbrace{\left(  A_{\bullet,\sim
v}\right)  _{\bullet,r}}_{=A_{\bullet,r}}\right)  =\det\left(  A_{\bullet,\sim
u}\mid A_{\bullet,r}\right)  =0
\]
(by Proposition \ref{prop.addcol.props3} \textbf{(c)}, applied to $q=u$).
Hence, Proposition \ref{prop.addcol.props3} \textbf{(e)} is proven in Case 1.

Let us now consider Case 2. In this case, we have $r\geq v$. Hence,
$r\in\left\{  v,v+1,\ldots,n-1\right\}  $ (since $r\in\left\{  1,2,\ldots
,n-1\right\}  $), so that $r+1\in\left\{  v+1,v+2,\ldots,n\right\}
\subseteq\left\{  1,2,\ldots,n\right\}  $. Furthermore, $r+1>r\geq v>u$ (since
$u<v$) and thus $r+1\neq u$. Hence, Proposition \ref{prop.addcol.props3}
\textbf{(c)} (applied to $r+1$ and $u$ instead of $r$ and $q$) yields
$\det\left(  A_{\bullet,\sim u}\mid A_{\bullet,r+1}\right)  =0$.

But recall that $r\in\left\{  v,v+1,\ldots,n-1\right\}  $. Thus, Proposition
\ref{prop.unrows.basics} \textbf{(e)} (applied to $m=n$ and $w=r$) yields
$\left(  A_{\bullet,\sim v}\right)  _{\bullet,r}=A_{\bullet,r+1}$. Hence,%
\[
\det\left(  A_{\bullet,\sim u}\mid\underbrace{\left(  A_{\bullet,\sim
v}\right)  _{\bullet,r}}_{=A_{\bullet,r+1}}\right)  =\det\left(
A_{\bullet,\sim u}\mid A_{\bullet,r+1}\right)  =0.
\]
Hence, Proposition \ref{prop.addcol.props3} \textbf{(e)} is proven in Case 2.

We have now proven Proposition \ref{prop.addcol.props3} \textbf{(e)} in each
of the two Cases 1 and 2. Hence, Proposition \ref{prop.addcol.props3}
\textbf{(e)} always holds.

\textbf{(f)} Let $u$ and $v$ be two elements of $\left\{  1,2,\ldots
,n\right\}  $ satisfying $u<v$. Then, $u<v$, so that $u\in\left\{
1,2,\ldots,v-1\right\}  $. Hence, Proposition \ref{prop.unrows.basics}
\textbf{(d)} (applied to $m=n$ and $w=u$) yields $\left(  A_{\bullet,\sim
v}\right)  _{\bullet,u}=A_{\bullet,u}$. Thus,%
\[
\det\left(  A_{\bullet,\sim u}\mid\underbrace{\left(  A_{\bullet,\sim
v}\right)  _{\bullet,u}}_{=A_{\bullet,u}}\right)  =\det\left(  A_{\bullet,\sim
u}\mid A_{\bullet,u}\right)  =\left(  -1\right)  ^{n+u}\det A
\]
(by Proposition \ref{prop.addcol.props3} \textbf{(b)}, applied to $q=u$).
Thus,%
\[
\left(  -1\right)  ^{u}\underbrace{\det\left(  A_{\bullet,\sim u}\mid\left(
A_{\bullet,\sim v}\right)  _{\bullet,u}\right)  }_{=\left(  -1\right)
^{n+u}\det A}=\underbrace{\left(  -1\right)  ^{u}\left(  -1\right)  ^{n+u}%
}_{\substack{=\left(  -1\right)  ^{u+\left(  n+u\right)  }=\left(  -1\right)
^{n}\\\text{(since }u+\left(  n+u\right)  =2u+n\equiv n\operatorname{mod}%
2\text{)}}}\det A=\left(  -1\right)  ^{n}\det A.
\]
This proves Proposition \ref{prop.addcol.props3} \textbf{(f)}.
\end{proof}
\end{vershort}

\begin{verlong}
In preparation for our solution to Exercise \ref{exe.prop.addcol.props}, we
shall prove a lemma:

\begin{lemma}
\label{lem.sol.prop.addcol.props.1}Let $n\in\mathbb{N}$ and $m\in\mathbb{N}$.
Let $A=\left(  a_{i,j}\right)  _{1\leq i\leq n,\ 1\leq j\leq m}\in
\mathbb{K}^{n\times m}$ be an $n\times m$-matrix. Let $v=\left(  v_{1}%
,v_{2},\ldots,v_{n}\right)  ^{T}\in\mathbb{K}^{n\times1}$ be a column vector
with $n$ entries.

\textbf{(a)} Every $q\in\left\{  1,2,\ldots,m\right\}  $ satisfies $\left(
A\mid v\right)  _{\bullet,q}=A_{\bullet,q}$.

\textbf{(b)} We have $\left(  A\mid v\right)  _{\bullet,m+1}=v$.

\textbf{(c)} Write the matrix $\left(  A\mid v\right)  $ in the form $\left(
A\mid v\right)  =\left(  b_{i,j}\right)  _{1\leq i\leq n,\ 1\leq j\leq m+1}$.
Then, we have%
\begin{equation}
b_{i,j}=a_{i,j} \label{eq.lem.sol.prop.addcol.props.1.1}%
\end{equation}
for every $i\in\left\{  1,2,\ldots,n\right\}  $ and every $j\in\left\{
1,2,\ldots,m\right\}  $. Moreover,%
\begin{equation}
b_{i,m+1}=v_{i} \label{eq.lem.sol.prop.addcol.props.1.2}%
\end{equation}
for every $i\in\left\{  1,2,\ldots,n\right\}  $.
\end{lemma}

\begin{noncompile}
[The following version of the proof of Lemma \ref{lem.sol.prop.addcol.props.1}
was written for the vershort environment, back when I thought that this lemma
would appear in this environment at all:]

\begin{proof}
[Proof of Lemma \ref{lem.sol.prop.addcol.props.1}.]The matrix $\left(  A\mid
v\right)  $ is defined as the $n\times\left(  m+1\right)  $-matrix whose $m+1$
columns are $A_{\bullet,1},A_{\bullet,2},\ldots,A_{\bullet,m},v$ (from left to
right). Thus, the first $m$ columns of the matrix $\left(  A\mid v\right)  $
are $A_{\bullet,1},A_{\bullet,2},\ldots,A_{\bullet,m}$, whereas the $\left(
m+1\right)  $-st column of $\left(  A\mid v\right)  $ is $v$.

\textbf{(a)} For every $q\in\left\{  1,2,\ldots,m\right\}  $, we have%
\begin{align*}
\left(  A\mid v\right)  _{\bullet,q}  &  =\left(  \text{the }q\text{-th column
of the matrix }\left(  A\mid v\right)  \right) \\
&  \ \ \ \ \ \ \ \ \ \ \left(  \text{by the definition of }\left(  A\mid
v\right)  _{\bullet,q}\right) \\
&  =A_{\bullet,q}\ \ \ \ \ \ \ \ \ \ \left(
\begin{array}
[c]{c}%
\text{since the first }m\text{ columns of the matrix }\left(  A\mid v\right)
\\
\text{are }A_{\bullet,1},A_{\bullet,2},\ldots,A_{\bullet,m}%
\end{array}
\right)  .
\end{align*}
This proves Lemma \ref{lem.sol.prop.addcol.props.1} \textbf{(a)}.

\textbf{(b)} The definition of $\left(  A\mid v\right)  _{\bullet,m+1}$ yields%
\[
\left(  A\mid v\right)  _{\bullet,m+1}=\left(  \text{the }\left(  m+1\right)
\text{-st column of the matrix }\left(  A\mid v\right)  \right)  =v
\]
(since the $\left(  m+1\right)  $-st column of $\left(  A\mid v\right)  $ is
$v$). This proves Lemma \ref{lem.sol.prop.addcol.props.1} \textbf{(b)}.

\textbf{(c)} [\textit{Proof of (\ref{eq.lem.sol.prop.addcol.props.1.1}):} We
have $\left(  A\mid v\right)  =\left(  b_{i,j}\right)  _{1\leq i\leq n,\ 1\leq
j\leq m+1}$. Thus, for every $i\in\left\{  1,2,\ldots,n\right\}  $ and
$j\in\left\{  1,2,\ldots,m\right\}  $, we have%
\begin{align*}
b_{i,j}  &  =\left(  \text{the }\left(  i,j\right)  \text{-th entry of the
matrix }\left(  A\mid v\right)  \right) \\
&  =\left(  \text{the }i\text{-th entry of }\underbrace{\text{the }j\text{-th
column of the matrix }\left(  A\mid v\right)  }_{\substack{=\left(  A\mid
v\right)  _{\bullet,j}=A_{\bullet,j}\\\text{(by Lemma
\ref{lem.sol.prop.addcol.props.1} \textbf{(a)},}\\\text{applied to
}q=j\text{)}}}\right) \\
&  =\left(  \text{the }i\text{-th entry of }\underbrace{A_{\bullet,j}%
}_{\substack{=\left(  \text{the }j\text{-th column of the matrix }A\right)
\\\text{(by the definition of }A_{\bullet,j}\text{)}}}\right) \\
&  =\left(  \text{the }i\text{-th entry of the }j\text{-th column of the
matrix }A\right) \\
&  =\left(  \text{the }\left(  i,j\right)  \text{-th entry of the matrix
}A\right)  =a_{i,j}.
\end{align*}
This proves (\ref{eq.lem.sol.prop.addcol.props.1.1}).]

[\textit{Proof of (\ref{eq.lem.sol.prop.addcol.props.1.2}):} We have $\left(
A\mid v\right)  =\left(  b_{i,j}\right)  _{1\leq i\leq n,\ 1\leq j\leq m+1}$.
Thus, for every $i\in\left\{  1,2,\ldots,n\right\}  $, we have%
\begin{align*}
b_{i,m+1}  &  =\left(  \text{the }\left(  i,m+1\right)  \text{-th entry of the
matrix }\left(  A\mid v\right)  \right) \\
&  =\left(  \text{the }i\text{-th entry of }\underbrace{\text{the }\left(
m+1\right)  \text{-th column of the matrix }\left(  A\mid v\right)
}_{\substack{=\left(  A\mid v\right)  _{\bullet,m+1}=v\\\text{(by Lemma
\ref{lem.sol.prop.addcol.props.1} \textbf{(b)})}}}\right) \\
&  =\left(  \text{the }i\text{-th entry of }v\right)  =v_{i}.
\end{align*}
This proves (\ref{eq.lem.sol.prop.addcol.props.1.2}).]

Thus, Lemma \ref{lem.sol.prop.addcol.props.1} \textbf{(c)} is proven.
\end{proof}
\end{noncompile}

\begin{proof}
[Proof of Lemma \ref{lem.sol.prop.addcol.props.1}.]The matrix $\left(  A\mid
v\right)  $ is defined as the $n\times\left(  m+1\right)  $-matrix whose $m+1$
columns are $A_{\bullet,1},A_{\bullet,2},\ldots,A_{\bullet,m},v$ (from left to right).

For each $q\in\left\{  1,2,\ldots,m+1\right\}  $, we have%
\begin{align}
&  \left(  \text{the }q\text{-th column of the matrix }\left(  A\mid v\right)
\right) \nonumber\\
&  =\left(  \text{the }q\text{-th entry of }\underbrace{\text{the list of the
columns of the matrix }\left(  A\mid v\right)  }_{\substack{=\left(
A_{\bullet,1},A_{\bullet,2},\ldots,A_{\bullet,m},v\right)  \\\text{(since
}\left(  A\mid v\right)  \text{ is the }n\times\left(  m+1\right)
\text{-matrix whose }m+1\\\text{columns are }A_{\bullet,1},A_{\bullet
,2},\ldots,A_{\bullet,m},v\text{ (from left to right))}}}\right) \nonumber\\
&  =\left(  \text{the }q\text{-th entry of }\left(  A_{\bullet,1}%
,A_{\bullet,2},\ldots,A_{\bullet,m},v\right)  \right) \nonumber\\
&  =%
\begin{cases}
A_{\bullet,q}, & \text{if }q<m+1;\\
v, & \text{if }q=m+1
\end{cases}
. \label{pf.lem.sol.prop.addcol.props.1.long.1}%
\end{align}

Now, for each $q\in\left\{  1,2,\ldots,m+1\right\}  $, we have%
\begin{align}
&  \left(  A\mid v\right)  _{\bullet,q}\nonumber\\
&  =\left(  \text{the }q\text{-th column of the matrix }\left(  A\mid
v\right)  \right) \nonumber\\
&  \ \ \ \ \ \ \ \ \ \ \left(  \text{by the definition of }\left(  A\mid
v\right)  _{\bullet,q}\right) \nonumber\\
&  =%
\begin{cases}
A_{\bullet,q}, & \text{if }q<m+1;\\
v, & \text{if }q=m+1
\end{cases}
\ \ \ \ \ \ \ \ \ \ \left(  \text{by
(\ref{pf.lem.sol.prop.addcol.props.1.long.1})}\right)  .
\label{pf.lem.sol.prop.addcol.props.1.long.2}%
\end{align}

\textbf{(a)} Let $q\in\left\{  1,2,\ldots,m\right\}  $. Then, $q\leq m<m+1$.
Moreover, $q\in\left\{  1,2,\ldots,m\right\}  \subseteq\left\{  1,2,\ldots
,m+1\right\}  $. Hence, (\ref{pf.lem.sol.prop.addcol.props.1.long.2}) yields%
\[
\left(  A\mid v\right)  _{\bullet,q}=%
\begin{cases}
A_{\bullet,q}, & \text{if }q<m+1;\\
v, & \text{if }q=m+1
\end{cases}
=A_{\bullet,q}\ \ \ \ \ \ \ \ \ \ \left(  \text{since }q<m+1\right)  .
\]
Thus, Lemma \ref{lem.sol.prop.addcol.props.1} \textbf{(a)} is proven.

\textbf{(b)} We have $\underbrace{m}_{\geq0}+1\geq1$, and thus $m+1\in\left\{
1,2,\ldots,m+1\right\}  $. Hence, (\ref{pf.lem.sol.prop.addcol.props.1.long.2}%
) (applied to $q=m+1$) yields%
\[
\left(  A\mid v\right)  _{\bullet,m+1}=%
\begin{cases}
A_{\bullet,m+1}, & \text{if }m+1<m+1;\\
v, & \text{if }m+1=m+1
\end{cases}
=v\ \ \ \ \ \ \ \ \ \ \left(  \text{since }m+1=m+1\right)  .
\]
Thus, Lemma \ref{lem.sol.prop.addcol.props.1} \textbf{(b)} is proven.

\textbf{(c)} We have $\left(  b_{i,j}\right)  _{1\leq i\leq n,\ 1\leq j\leq
m+1}=\left(  A\mid v\right)  $. Thus, for every $i\in\left\{  1,2,\ldots
,n\right\}  $ and $j\in\left\{  1,2,\ldots,m+1\right\}  $, we have%
\begin{align}
b_{i,j}  &  =\left(  \text{the }\left(  i,j\right)  \text{-th entry of the
matrix }\left(  A\mid v\right)  \right)
\label{pf.lem.sol.prop.addcol.props.1.long.c.0}\\
&  =\left(  \text{the }i\text{-th entry of }\underbrace{\text{the }j\text{-th
column of the matrix }\left(  A\mid v\right)  }_{\substack{=\left(  A\mid
v\right)  _{\bullet,j}\\\text{(since }\left(  A\mid v\right)  _{\bullet
,j}\text{ is the }j\text{-th column of the matrix }\left(  A\mid v\right)
\\\text{(by the definition of }\left(  A\mid v\right)  _{\bullet,j}\text{))}%
}}\right) \nonumber\\
&  =\left(  \text{the }i\text{-th entry of }\left(  A\mid v\right)
_{\bullet,j}\right)  . \label{pf.lem.sol.prop.addcol.props.1.long.c.0b}%
\end{align}

[\textit{Proof of (\ref{eq.lem.sol.prop.addcol.props.1.1}):} We have $\left(
a_{i,j}\right)  _{1\leq i\leq n,\ 1\leq j\leq m}=A$. Hence, for every
$i\in\left\{  1,2,\ldots,n\right\}  $ and $j\in\left\{  1,2,\ldots,m\right\}
$, we have%
\begin{equation}
a_{i,j}=\left(  \text{the }\left(  i,j\right)  \text{-th entry of the matrix
}A\right)  . \label{pf.lem.sol.prop.addcol.props.1.long.c.1}%
\end{equation}

Now, let $i\in\left\{  1,2,\ldots,n\right\}  $ and $j\in\left\{
1,2,\ldots,m\right\}  $. Then, $j\in\left\{  1,2,\ldots,m\right\}
\subseteq\left\{  1,2,\ldots,m+1\right\}  $. Thus,
(\ref{pf.lem.sol.prop.addcol.props.1.long.c.0b}) yields%
\begin{align*}
b_{i,j}  &  =\left(  \text{the }i\text{-th entry of }\underbrace{\left(  A\mid
v\right)  _{\bullet,j}}_{\substack{=A_{\bullet,j}\\\text{(by Lemma
\ref{lem.sol.prop.addcol.props.1} \textbf{(a)},}\\\text{applied to
}q=j\text{)}}}\right) \\
&  =\left(  \text{the }i\text{-th entry of }\underbrace{A_{\bullet,j}%
}_{\substack{=\left(  \text{the }j\text{-th column of the matrix }A\right)
\\\text{(since }A_{\bullet,j}\text{ is defined as the }j\text{-th column of
the matrix }A\text{)}}}\right) \\
&  =\left(  \text{the }i\text{-th entry of the }j\text{-th column of the
matrix }A\right) \\
&  =\left(  \text{the }\left(  i,j\right)  \text{-th entry of the matrix
}A\right)  =a_{i,j}\ \ \ \ \ \ \ \ \ \ \left(  \text{by
(\ref{pf.lem.sol.prop.addcol.props.1.long.c.1})}\right)  .
\end{align*}
This proves (\ref{eq.lem.sol.prop.addcol.props.1.1}).]

[\textit{Proof of (\ref{eq.lem.sol.prop.addcol.props.1.2}):} We have
$v=\left(  v_{1},v_{2},\ldots,v_{n}\right)  ^{T}=\left(
\begin{array}
[c]{c}%
v_{1}\\
v_{2}\\
\vdots\\
v_{n}%
\end{array}
\right)  $. Hence, for every $i\in\left\{  1,2,\ldots,n\right\}  $, we have%
\begin{equation}
\left(  \text{the }i\text{-th entry of }v\right)  =v_{i}.
\label{pf.lem.sol.prop.addcol.props.1.long.c.2}%
\end{equation}

Now, let $i\in\left\{  1,2,\ldots,n\right\}  $. We have $\underbrace{m}%
_{\geq0}+1\geq1$, and thus $m+1\in\left\{  1,2,\ldots,m+1\right\}  $. Thus,
(\ref{pf.lem.sol.prop.addcol.props.1.long.c.0b}) (applied to $j=m+1$) yields%
\begin{align*}
b_{i,m+1}  &  =\left(  \text{the }i\text{-th entry of }\underbrace{\left(
A\mid v\right)  _{\bullet,m+1}}_{\substack{=v\\\text{(by Lemma
\ref{lem.sol.prop.addcol.props.1} \textbf{(b)})}}}\right) \\
&  =\left(  \text{the }i\text{-th entry of }v\right)  =v_{i}%
\ \ \ \ \ \ \ \ \ \ \left(  \text{by
(\ref{pf.lem.sol.prop.addcol.props.1.long.c.2})}\right)  .
\end{align*}
This proves (\ref{eq.lem.sol.prop.addcol.props.1.2}).]

Now, we have proven both (\ref{eq.lem.sol.prop.addcol.props.1.1}) and
(\ref{eq.lem.sol.prop.addcol.props.1.2}). Thus, Lemma
\ref{lem.sol.prop.addcol.props.1} \textbf{(c)} is proven.
\end{proof}

\begin{proof}
[Proof of Proposition \ref{prop.addcol.props1}.]Write the $n\times m$-matrix
$A$ in the form $A=\left(  a_{i,j}\right)  _{1\leq i\leq n,\ 1\leq j\leq m}%
\in\mathbb{K}^{n\times m}$. Write the column vector $v\in\mathbb{K}^{n\times
1}$ in the form $v=\left(
\begin{array}
[c]{c}%
v_{1}\\
v_{2}\\
\vdots\\
v_{n}%
\end{array}
\right)  $. Thus, $v=\left(
\begin{array}
[c]{c}%
v_{1}\\
v_{2}\\
\vdots\\
v_{n}%
\end{array}
\right)  =\left(  v_{1},v_{2},\ldots,v_{n}\right)  ^{T}$.

We shall use the notation introduced in Definition \ref{def.sol.unrows.d}
\textbf{(a)}.

\textbf{(a)} Lemma \ref{lem.sol.prop.addcol.props.1} \textbf{(a)} shows that
every $q\in\left\{  1,2,\ldots,m\right\}  $ satisfies $\left(  A\mid v\right)
_{\bullet,q}=A_{\bullet,q}$. This proves Proposition \ref{prop.addcol.props1}
\textbf{(a)}.

\textbf{(b)} Lemma \ref{lem.sol.prop.addcol.props.1} \textbf{(b)} shows that
$\left(  A\mid v\right)  _{\bullet,m+1}=v$. This proves Proposition
\ref{prop.addcol.props1} \textbf{(b)}.

\textbf{(c)} Write the matrix $\left(  A\mid v\right)  $ in the form $\left(
A\mid v\right)  =\left(  b_{i,j}\right)  _{1\leq i\leq n,\ 1\leq j\leq m+1}$.
(This is possible since $\left(  A\mid v\right)  $ is an $n\times\left(
m+1\right)  $-matrix.) Then, we have%
\begin{equation}
b_{i,j}=a_{i,j} \label{pf.prop.addcol.props1.c.1}%
\end{equation}
for every $i\in\left\{  1,2,\ldots,n\right\}  $ and every $j\in\left\{
1,2,\ldots,m\right\}  $ (by (\ref{eq.lem.sol.prop.addcol.props.1.1})).
Moreover,%
\begin{equation}
b_{i,m+1}=v_{i} \label{pf.prop.addcol.props1.c.2}%
\end{equation}
for every $i\in\left\{  1,2,\ldots,n\right\}  $ (by
(\ref{eq.lem.sol.prop.addcol.props.1.2})).

Let $q\in\left\{  1,2,\ldots,m\right\}  $. The matrix $A_{\bullet,\sim q}$ is
an $n\times\left(  m-1\right)  $-matrix (since $A$ is an $n\times m$-matrix).
Hence, $\left(  A_{\bullet,\sim q}\mid v\right)  $ is an $n\times\left(
\left(  m-1\right)  +1\right)  $-matrix. In other words, $\left(
A_{\bullet,\sim q}\mid v\right)  $ is an $n\times m$-matrix (since $\left(
m-1\right)  +1=m$).

Proposition \ref{prop.sol.unirows.d} \textbf{(b)} (applied to $v=q$) shows
that $A_{\bullet,\sim q}=\left(  a_{i,\mathbf{d}_{q}\left(  j\right)
}\right)  _{1\leq i\leq n,\ 1\leq j\leq m-1}$. But recall that $\left(  A\mid
v\right)  =\left(  b_{i,j}\right)  _{1\leq i\leq n,\ 1\leq j\leq m+1}$. Hence,
Proposition \ref{prop.sol.unirows.d} \textbf{(b)} (applied to $m+1$, $\left(
A\mid v\right)  $, $b_{i,j}$ and $q$ instead of $m$, $A$, $a_{i,j}$ and $v$)
shows that
\begin{equation}
\left(  A\mid v\right)  _{\bullet,\sim q}=\left(  b_{i,\mathbf{d}_{q}\left(
j\right)  }\right)  _{1\leq i\leq n,\ 1\leq j\leq\left(  m+1\right)
-1}=\left(  b_{i,\mathbf{d}_{q}\left(  j\right)  }\right)  _{1\leq i\leq
n,\ 1\leq j\leq m} \label{pf.prop.addcol.props1.c.mat1}%
\end{equation}
(since $\left(  m+1\right)  -1=m$).

On the other hand, write the matrix $\left(  A_{\bullet,\sim q}\mid v\right)
$ in the form
\[
\left(  A_{\bullet,\sim q}\mid v\right)  =\left(  c_{i,j}\right)  _{1\leq
i\leq n,\ 1\leq j\leq\left(  m-1\right)  +1}.
\]
(This is possible since $\left(  A_{\bullet,\sim q}\mid v\right)  $ is an
$n\times\left(  \left(  m-1\right)  +1\right)  $-matrix.) Recall that
$A_{\bullet,\sim q}=\left(  a_{i,\mathbf{d}_{q}\left(  j\right)  }\right)
_{1\leq i\leq n,\ 1\leq j\leq m-1}$. Thus, we have%
\begin{equation}
c_{i,j}=a_{i,\mathbf{d}_{q}\left(  j\right)  }
\label{pf.prop.addcol.props1.c.3}%
\end{equation}
for every $i\in\left\{  1,2,\ldots,n\right\}  $ and every $j\in\left\{
1,2,\ldots,m-1\right\}  $ (by (\ref{eq.lem.sol.prop.addcol.props.1.1}),
applied to $m-1$, $A_{\bullet,\sim q}$, $a_{i,\mathbf{d}_{q}\left(  j\right)
}$ and $c_{i,j}$ instead of $m$, $A$, $a_{i,j}$ and $b_{i,j}$). Moreover,%
\begin{equation}
c_{i,\left(  m-1\right)  +1}=v_{i} \label{pf.prop.addcol.props1.c.4}%
\end{equation}
for every $i\in\left\{  1,2,\ldots,n\right\}  $ (by
(\ref{eq.lem.sol.prop.addcol.props.1.2}), applied to $m-1$, $A_{\bullet,\sim
q}$, $a_{i,\mathbf{d}_{q}\left(  j\right)  }$ and $c_{i,j}$ instead of $m$,
$A$, $a_{i,j}$ and $b_{i,j}$).

Now,
\begin{equation}
\left(  A_{\bullet,\sim q}\mid v\right)  =\left(  c_{i,j}\right)  _{1\leq
i\leq n,\ 1\leq j\leq\left(  m-1\right)  +1}=\left(  c_{i,j}\right)  _{1\leq
i\leq n,\ 1\leq j\leq m} \label{pf.prop.addcol.props1.c.mat2}%
\end{equation}
(since $\left(  m-1\right)  +1=m$). We shall now show that%
\begin{equation}
b_{i,\mathbf{d}_{q}\left(  j\right)  }=c_{i,j}
\label{pf.prop.addcol.props1.c.10}%
\end{equation}
for every $i\in\left\{  1,2,\ldots,n\right\}  $ and $j\in\left\{
1,2,\ldots,m\right\}  $.

[\textit{Proof of (\ref{pf.prop.addcol.props1.c.10}):} Let $i\in\left\{
1,2,\ldots,n\right\}  $ and $j\in\left\{  1,2,\ldots,m\right\}  $. We must
prove the equality (\ref{pf.prop.addcol.props1.c.10}). We are in one of the
following two cases:

\textit{Case 1:} We have $j\neq m$.

\textit{Case 2:} We have $j=m$.

Let us first consider Case 1. In this case, we have $j\neq m$. Combined with
$j\in\left\{  1,2,\ldots,m\right\}  $, this yields $j\in\left\{
1,2,\ldots,m\right\}  \setminus\left\{  m\right\}  =\left\{  1,2,\ldots
,m-1\right\}  $. Hence, $j+1\in\left\{  2,3,\ldots,m\right\}  $, so that
$j+1\leq m$. Now, the definition of $\mathbf{d}_{q}\left(  j\right)  $ yields%
\begin{align*}
\mathbf{d}_{q}\left(  j\right)   &  =%
\begin{cases}
j, & \text{if }j<q;\\
j+1, & \text{if }j\geq q
\end{cases}
\\
&  \leq%
\begin{cases}
j+1, & \text{if }j<q;\\
j+1, & \text{if }j\geq q
\end{cases}
\ \ \ \ \ \ \ \ \ \ \left(
\begin{array}
[c]{c}%
\text{since }j\leq j+1\text{ in the case when }j<q\text{,}\\
\text{and since }j+1\leq j+1\text{ in the case when }j\geq q
\end{array}
\right) \\
&  =j+1\leq m.
\end{align*}
Combining this with%
\begin{align*}
\mathbf{d}_{q}\left(  j\right)   &  =%
\begin{cases}
j, & \text{if }j<q;\\
j+1, & \text{if }j\geq q
\end{cases}
\\
&  \geq%
\begin{cases}
j, & \text{if }j<q;\\
j, & \text{if }j\geq q
\end{cases}
\ \ \ \ \ \ \ \ \ \ \left(
\begin{array}
[c]{c}%
\text{since }j\geq j\text{ in the case when }j<q\text{,}\\
\text{and since }j+1\geq j\text{ in the case when }j\geq q
\end{array}
\right) \\
&  =j\geq1\ \ \ \ \ \ \ \ \ \ \left(  \text{since }j\in\left\{  1,2,\ldots
,m\right\}  \right)  ,
\end{align*}
we obtain $\mathbf{d}_{q}\left(  j\right)  \in\left\{  1,2,\ldots,m\right\}
$. Hence, (\ref{pf.prop.addcol.props1.c.1}) (applied to $\mathbf{d}_{q}\left(
j\right)  $ instead of $j$) shows that $b_{i,\mathbf{d}_{q}\left(  j\right)
}=a_{i,\mathbf{d}_{q}\left(  j\right)  }$. But
(\ref{pf.prop.addcol.props1.c.3}) yields $c_{i,j}=a_{i,\mathbf{d}_{q}\left(
j\right)  }$. Thus, $b_{i,\mathbf{d}_{q}\left(  j\right)  }=a_{i,\mathbf{d}%
_{q}\left(  j\right)  }=c_{i,j}$. Hence, (\ref{pf.prop.addcol.props1.c.10}) is
proven in Case 1.

Let us now consider Case 2. In this case, we have $j=m$. Hence,
\begin{align*}
\mathbf{d}_{q}\left(  \underbrace{j}_{=m}\right)   &  =\mathbf{d}_{q}\left(
m\right)  =%
\begin{cases}
m, & \text{if }m<q;\\
m+1, & \text{if }m\geq q
\end{cases}
\ \ \ \ \ \ \ \ \ \ \left(  \text{by the definition of }\mathbf{d}_{q}\left(
m\right)  \right) \\
&  =m+1\ \ \ \ \ \ \ \ \ \ \left(  \text{since }m\geq q\text{ (since }q\leq
m\text{ (because }q\in\left\{  1,2,\ldots,m\right\}  \text{))}\right)  .
\end{align*}
Hence, $b_{i,\mathbf{d}_{q}\left(  j\right)  }=b_{i,m+1}=v_{i}$ (by
(\ref{pf.prop.addcol.props1.c.2})). But from $j=m=\left(  m-1\right)  +1$, we
obtain $c_{i,j}=c_{i,\left(  m-1\right)  +1}=v_{i}$ (by
(\ref{pf.prop.addcol.props1.c.4})). Thus, $b_{i,\mathbf{d}_{q}\left(
j\right)  }=v_{i}=c_{i,j}$. Hence, (\ref{pf.prop.addcol.props1.c.10}) is
proven in Case 2.

We have now proven the equality (\ref{pf.prop.addcol.props1.c.10}) in each of
the two Cases 1 and 2. Since these two Cases cover all possibilities, this
shows that (\ref{pf.prop.addcol.props1.c.10}) always holds. Thus,
(\ref{pf.prop.addcol.props1.c.10}) is proven.]

Now, (\ref{pf.prop.addcol.props1.c.mat1}) becomes%
\[
\left(  A\mid v\right)  _{\bullet,\sim q}=\left(  \underbrace{b_{i,\mathbf{d}%
_{q}\left(  j\right)  }}_{\substack{=c_{i,j}\\\text{(by
(\ref{pf.prop.addcol.props1.c.10}))}}}\right)  _{1\leq i\leq n,\ 1\leq j\leq
m}=\left(  c_{i,j}\right)  _{1\leq i\leq n,\ 1\leq j\leq m}=\left(
A_{\bullet,\sim q}\mid v\right)
\]
(by (\ref{pf.prop.addcol.props1.c.mat2})). This proves Proposition
\ref{prop.addcol.props1} \textbf{(c)}.

\textbf{(d)} Write the matrix $\left(  A\mid v\right)  $ in the form $\left(
A\mid v\right)  =\left(  b_{i,j}\right)  _{1\leq i\leq n,\ 1\leq j\leq m+1}$.
(This is possible since $\left(  A\mid v\right)  $ is an $n\times\left(
m+1\right)  $-matrix.) Then, we have%
\begin{equation}
b_{i,j}=a_{i,j} \label{pf.prop.addcol.props1.d.1}%
\end{equation}
for every $i\in\left\{  1,2,\ldots,n\right\}  $ and every $j\in\left\{
1,2,\ldots,m\right\}  $ (by (\ref{eq.lem.sol.prop.addcol.props.1.1}))

We have $m+1\in\left\{  1,2,\ldots,m+1\right\}  $ (since $\underbrace{m}%
_{\geq0}+1\geq1$). Thus, the matrix $\left(  A\mid v\right)  _{\bullet
,\sim\left(  m+1\right)  }$ is well-defined (since $\left(  A\mid v\right)  $
is an $n\times\left(  m+1\right)  $-matrix). This matrix $\left(  A\mid
v\right)  _{\bullet,\sim\left(  m+1\right)  }$ is an $n\times\left(  \left(
m+1\right)  -1\right)  $-matrix (since $\left(  A\mid v\right)  $ is an
$n\times\left(  m+1\right)  $-matrix). In other words, $\left(  A\mid
v\right)  _{\bullet,\sim\left(  m+1\right)  }$ is an $n\times m$-matrix (since
$\left(  m+1\right)  -1=m$). Proposition \ref{prop.sol.unirows.d} \textbf{(b)}
(applied to $m+1$, $\left(  A\mid v\right)  $, $b_{i,j}$ and $m+1$ instead of
$m$, $A$, $a_{i,j}$ and $v$) shows that
\begin{equation}
\left(  A\mid v\right)  _{\bullet,\sim\left(  m+1\right)  }=\left(
b_{i,\mathbf{d}_{m+1}\left(  j\right)  }\right)  _{1\leq i\leq n,\ 1\leq
j\leq\left(  m+1\right)  -1}=\left(  b_{i,\mathbf{d}_{m+1}\left(  j\right)
}\right)  _{1\leq i\leq n,\ 1\leq j\leq m} \label{pf.prop.addcol.props1.d.2}%
\end{equation}
(since $\left(  m+1\right)  -1=m$).

But every $i\in\left\{  1,2,\ldots,n\right\}  $ and $j\in\left\{
1,2,\ldots,m\right\}  $ satisfy
\begin{equation}
b_{i,\mathbf{d}_{m+1}\left(  j\right)  }=a_{i,j}
\label{pf.prop.addcol.props1.d.3}%
\end{equation}
\footnote{\textit{Proof of (\ref{pf.prop.addcol.props1.d.3}):} Let
$i\in\left\{  1,2,\ldots,n\right\}  $ and $j\in\left\{  1,2,\ldots,m\right\}
$. From $j\in\left\{  1,2,\ldots,m\right\}  $, we obtain $j\leq m<m+1$. Now,
the definition of $\mathbf{d}_{m+1}\left(  j\right)  $ yields $\mathbf{d}%
_{m+1}\left(  j\right)  =%
\begin{cases}
j, & \text{if }j<m+1;\\
j+1, & \text{if }j\geq m+1
\end{cases}
=j$ (since $j<m+1$). Thus, $b_{i,\mathbf{d}_{m+1}\left(  j\right)  }%
=b_{i,j}=a_{i,j}$ (by (\ref{pf.prop.addcol.props1.d.1})). This proves
(\ref{pf.prop.addcol.props1.d.3}).}. Hence, $\left(  b_{i,\mathbf{d}%
_{m+1}\left(  j\right)  }\right)  _{1\leq i\leq n,\ 1\leq j\leq m}=\left(
a_{i,j}\right)  _{1\leq i\leq n,\ 1\leq j\leq m}=A$. Thus,
(\ref{pf.prop.addcol.props1.d.2}) becomes%
\[
\left(  A\mid v\right)  _{\bullet,\sim\left(  m+1\right)  }=\left(
b_{i,\mathbf{d}_{m+1}\left(  j\right)  }\right)  _{1\leq i\leq n,\ 1\leq j\leq
m}=A.
\]
This proves Proposition \ref{prop.addcol.props1} \textbf{(d)}.

\textbf{(e)} Write the matrix $\left(  A\mid v\right)  $ in the form $\left(
A\mid v\right)  =\left(  b_{i,j}\right)  _{1\leq i\leq n,\ 1\leq j\leq m+1}$.
(This is possible since $\left(  A\mid v\right)  $ is an $n\times\left(
m+1\right)  $-matrix.) Then, we have%
\begin{equation}
b_{i,j}=a_{i,j} \label{pf.prop.addcol.props1.e.1}%
\end{equation}
for every $i\in\left\{  1,2,\ldots,n\right\}  $ and every $j\in\left\{
1,2,\ldots,m\right\}  $ (by (\ref{eq.lem.sol.prop.addcol.props.1.1})).
Moreover,%
\begin{equation}
b_{i,m+1}=v_{i} \label{pf.prop.addcol.props1.e.2}%
\end{equation}
for every $i\in\left\{  1,2,\ldots,n\right\}  $ (by
(\ref{eq.lem.sol.prop.addcol.props.1.2})).

Let $p\in\left\{  1,2,\ldots,n\right\}  $. The matrix $A_{\sim p,\bullet}$ is
an $\left(  n-1\right)  \times m$-matrix (since $A$ is an $n\times m$-matrix).
Also, the matrix $v_{\sim p,\bullet}$ is an $\left(  n-1\right)  \times
1$-matrix (since $v$ is an $n\times1$-matrix). Hence, $\left(  A_{\sim
p,\bullet}\mid v_{\sim p,\bullet}\right)  $ is an $\left(  n-1\right)
\times\left(  m+1\right)  $-matrix.

Proposition \ref{prop.sol.unirows.d} \textbf{(a)} (applied to $u=p$) shows
that $A_{\sim p,\bullet}=\left(  a_{\mathbf{d}_{p}\left(  i\right)
,j}\right)  _{1\leq i\leq n-1,\ 1\leq j\leq m}$. But recall that $\left(
A\mid v\right)  =\left(  b_{i,j}\right)  _{1\leq i\leq n,\ 1\leq j\leq m+1}$.
Hence, Proposition \ref{prop.sol.unirows.d} \textbf{(a)} (applied to $m+1$,
$\left(  A\mid v\right)  $, $b_{i,j}$ and $p$ instead of $m$, $A$, $a_{i,j}$
and $u$) shows that
\begin{equation}
\left(  A\mid v\right)  _{\sim p,\bullet}=\left(  b_{\mathbf{d}_{p}\left(
i\right)  ,j}\right)  _{1\leq i\leq n-1,\ 1\leq j\leq m+1}.
\label{pf.prop.addcol.props1.e.3}%
\end{equation}

Also, $v=\left(
\begin{array}
[c]{c}%
v_{1}\\
v_{2}\\
\vdots\\
v_{n}%
\end{array}
\right)  =\left(  v_{i}\right)  _{1\leq i\leq n,\ 1\leq j\leq1}$. Hence,
Proposition \ref{prop.sol.unirows.d} \textbf{(a)} (applied to $1$, $v$,
$v_{i}$ and $p$ instead of $m$, $A$, $a_{i,j}$ and $u$) yields
\begin{equation}
v_{\sim p,\bullet}=\left(  v_{\mathbf{d}_{p}\left(  i\right)  }\right)
_{1\leq i\leq n-1,\ 1\leq j\leq1}=\left(
\begin{array}
[c]{c}%
v_{\mathbf{d}_{p}\left(  1\right)  }\\
v_{\mathbf{d}_{p}\left(  2\right)  }\\
\vdots\\
v_{\mathbf{d}_{p}\left(  n-1\right)  }%
\end{array}
\right)  =\left(  v_{\mathbf{d}_{p}\left(  1\right)  },v_{\mathbf{d}%
_{p}\left(  2\right)  },\ldots,v_{\mathbf{d}_{p}\left(  n-1\right)  }\right)
^{T}. \label{pf.prop.addcol.props1.e.4}%
\end{equation}

On the other hand, write the matrix $\left(  A_{\sim p,\bullet}\mid v_{\sim
p,\bullet}\right)  $ in the form
\begin{equation}
\left(  A_{\sim p,\bullet}\mid v_{\sim p,\bullet}\right)  =\left(
c_{i,j}\right)  _{1\leq i\leq n-1,\ 1\leq j\leq m+1}.
\label{pf.prop.addcol.props1.e.mat2}%
\end{equation}
(This is possible since $\left(  A_{\sim p,\bullet}\mid v_{\sim p,\bullet
}\right)  $ is an $\left(  n-1\right)  \times\left(  m+1\right)  $-matrix.)
Recall that $A_{\sim p,\bullet}=\left(  a_{\mathbf{d}_{p}\left(  i\right)
,j}\right)  _{1\leq i\leq n-1,\ 1\leq j\leq m}$ and $v_{\sim p,\bullet
}=\left(  v_{\mathbf{d}_{p}\left(  1\right)  },v_{\mathbf{d}_{p}\left(
2\right)  },\ldots,v_{\mathbf{d}_{p}\left(  n-1\right)  }\right)  ^{T}$. Thus,
we have%
\begin{equation}
c_{i,j}=a_{\mathbf{d}_{p}\left(  i\right)  ,j}
\label{pf.prop.addcol.props1.e.5}%
\end{equation}
for every $i\in\left\{  1,2,\ldots,n-1\right\}  $ and every $j\in\left\{
1,2,\ldots,m\right\}  $ (by (\ref{eq.lem.sol.prop.addcol.props.1.1}), applied
to $n-1$, $A_{\sim p,\bullet}$, $a_{\mathbf{d}_{p}\left(  i\right)  ,j}$,
$v_{\sim p,\bullet}$, $v_{\mathbf{d}_{p}\left(  i\right)  }$ and $c_{i,j}$
instead of $n$, $A$, $a_{i,j}$, $v$, $v_{i}$ and $b_{i,j}$). Moreover,%
\begin{equation}
c_{i,m+1}=v_{\mathbf{d}_{p}\left(  i\right)  }
\label{pf.prop.addcol.props1.e.6}%
\end{equation}
for every $i\in\left\{  1,2,\ldots,n-1\right\}  $ (by
(\ref{eq.lem.sol.prop.addcol.props.1.2}), applied to $n-1$, $A_{\sim
p,\bullet}$, $a_{\mathbf{d}_{p}\left(  i\right)  ,j}$, $v_{\sim p,\bullet}$,
$v_{\mathbf{d}_{p}\left(  i\right)  }$ and $c_{i,j}$ instead of $n$, $A$,
$a_{i,j}$, $v$, $v_{i}$ and $b_{i,j}$).

We shall now show that%
\begin{equation}
b_{\mathbf{d}_{p}\left(  i\right)  ,j}=c_{i,j}
\label{pf.prop.addcol.props1.e.10}%
\end{equation}
for every $i\in\left\{  1,2,\ldots,n-1\right\}  $ and $j\in\left\{
1,2,\ldots,m+1\right\}  $.

[\textit{Proof of (\ref{pf.prop.addcol.props1.e.10}):} Let $i\in\left\{
1,2,\ldots,n-1\right\}  $ and $j\in\left\{  1,2,\ldots,m+1\right\}  $. We must
prove the equality (\ref{pf.prop.addcol.props1.e.10}). We have $\mathbf{d}%
_{p}\left(  i\right)  \in\left\{  1,2,\ldots,n\right\}  $%
\ \ \ \ \footnote{\textit{Proof.} From $i\in\left\{  1,2,\ldots,n-1\right\}
$, we obtain $i\leq n-1$. Now, the definition of $\mathbf{d}_{p}\left(
i\right)  $ yields%
\begin{align*}
\mathbf{d}_{p}\left(  i\right)   &  =%
\begin{cases}
i, & \text{if }i<p;\\
i+1, & \text{if }i\geq p
\end{cases}
\\
&  \leq%
\begin{cases}
i+1, & \text{if }i<p;\\
i+1, & \text{if }i\geq p
\end{cases}
\ \ \ \ \ \ \ \ \ \ \left(
\begin{array}
[c]{c}%
\text{since }i\leq i+1\text{ in the case when }i<p\text{,}\\
\text{and since }i+1\leq i+1\text{ in the case when }i\geq p
\end{array}
\right) \\
&  =i+1\leq n\ \ \ \ \ \ \ \ \ \ \left(  \text{since }i\leq n-1\right)  .
\end{align*}
Combined with%
\begin{align*}
\mathbf{d}_{p}\left(  i\right)   &  =%
\begin{cases}
i, & \text{if }i<p;\\
i+1, & \text{if }i\geq p
\end{cases}
\\
&  \geq%
\begin{cases}
i, & \text{if }i<p;\\
i, & \text{if }i\geq p
\end{cases}
\ \ \ \ \ \ \ \ \ \ \left(
\begin{array}
[c]{c}%
\text{since }i\geq i\text{ in the case when }i<p\text{,}\\
\text{and since }i+1\geq i\text{ in the case when }i\geq p
\end{array}
\right) \\
&  =i\geq1\ \ \ \ \ \ \ \ \ \ \left(  \text{since }i\in\left\{  1,2,\ldots
,n-1\right\}  \right)  ,
\end{align*}
this yields $\mathbf{d}_{p}\left(  i\right)  \in\left\{  1,2,\ldots,n\right\}
$. Qed.}.

We are in one of the following two cases:

\textit{Case 1:} We have $j\neq m+1$.

\textit{Case 2:} We have $j=m+1$.

Let us first consider Case 1. In this case, we have $j\neq m+1$. Combined with
$j\in\left\{  1,2,\ldots,m+1\right\}  $, this yields $j\in\left\{
1,2,\ldots,m+1\right\}  \setminus\left\{  m+1\right\}  =\left\{
1,2,\ldots,m\right\}  $. Hence, (\ref{pf.prop.addcol.props1.e.5}) yields
$c_{i,j}=a_{\mathbf{d}_{p}\left(  i\right)  ,j}$. Thus, $a_{\mathbf{d}%
_{p}\left(  i\right)  ,j}=c_{i,j}$. But (\ref{pf.prop.addcol.props1.e.1})
(applied to $\mathbf{d}_{p}\left(  i\right)  $ instead of $i$) yields
$b_{\mathbf{d}_{p}\left(  i\right)  ,j}=a_{\mathbf{d}_{p}\left(  i\right)
,j}$ (since $\mathbf{d}_{p}\left(  i\right)  \in\left\{  1,2,\ldots,n\right\}
$). Hence, $b_{\mathbf{d}_{p}\left(  i\right)  ,j}=a_{\mathbf{d}_{p}\left(
i\right)  ,j}=c_{i,j}$. Hence, (\ref{pf.prop.addcol.props1.e.10}) is proven in
Case 1.

Let us now consider Case 2. In this case, we have $j=m+1$. Hence,
$c_{i,j}=c_{i,m+1}=v_{\mathbf{d}_{p}\left(  i\right)  }$ (by
(\ref{pf.prop.addcol.props1.e.6})). Thus, $v_{\mathbf{d}_{p}\left(  i\right)
}=c_{i,j}$. But (\ref{pf.prop.addcol.props1.e.2}) (applied to $\mathbf{d}%
_{p}\left(  i\right)  $ instead of $i$) yields $b_{\mathbf{d}_{p}\left(
i\right)  ,m+1}=v_{\mathbf{d}_{p}\left(  i\right)  }$ (since $\mathbf{d}%
_{p}\left(  i\right)  \in\left\{  1,2,\ldots,n\right\}  $). Now, from $j=m+1$,
we obtain $b_{\mathbf{d}_{p}\left(  i\right)  ,j}=b_{\mathbf{d}_{p}\left(
i\right)  ,m+1}=v_{\mathbf{d}_{p}\left(  i\right)  }=c_{i,j}$. Hence,
(\ref{pf.prop.addcol.props1.e.10}) is proven in Case 2.

We have now proven the equality (\ref{pf.prop.addcol.props1.e.10}) in each of
the two Cases 1 and 2. Since these two Cases cover all possibilities, this
shows that (\ref{pf.prop.addcol.props1.e.10}) always holds. Thus,
(\ref{pf.prop.addcol.props1.e.10}) is proven.]

Now, (\ref{pf.prop.addcol.props1.e.3}) becomes%
\[
\left(  A\mid v\right)  _{\sim p,\bullet}=\left(  \underbrace{b_{\mathbf{d}%
_{p}\left(  i\right)  ,j}}_{\substack{=c_{i,j}\\\text{(by
(\ref{pf.prop.addcol.props1.e.10}))}}}\right)  _{1\leq i\leq n-1,\ 1\leq j\leq
m+1}=\left(  c_{i,j}\right)  _{1\leq i\leq n-1,\ 1\leq j\leq m+1}=\left(
A_{\sim p,\bullet}\mid v_{\sim p,\bullet}\right)
\]
(by (\ref{pf.prop.addcol.props1.e.mat2})). This proves Proposition
\ref{prop.addcol.props1} \textbf{(e)}.

\textbf{(f)} Let $p\in\left\{  1,2,\ldots,n\right\}  $. Clearly,
$m+1\in\left\{  1,2,\ldots,m+1\right\}  $. Also, $\left(  A\mid v\right)  $ is
an $n\times\left(  m+1\right)  $-matrix. Proposition \ref{prop.unrows.basics}
\textbf{(c)} (applied to $m+1$, $\left(  A\mid v\right)  $, $p$ and $m+1$
instead of $m$, $A$, $u$ and $v$) yields
\[
\left(  \left(  A\mid v\right)  _{\bullet,\sim\left(  m+1\right)  }\right)
_{\sim p,\bullet}=\left(  \left(  A\mid v\right)  _{\sim p,\bullet}\right)
_{\bullet,\sim\left(  m+1\right)  }=\left(  A\mid v\right)  _{\sim
p,\sim\left(  m+1\right)  }.
\]
Hence,%
\[
\left(  A\mid v\right)  _{\sim p,\sim\left(  m+1\right)  }=\left(
\underbrace{\left(  A\mid v\right)  _{\bullet,\sim\left(  m+1\right)  }%
}_{\substack{=A\\\text{(by Proposition \ref{prop.addcol.props1} \textbf{(d)}%
)}}}\right)  _{\sim p,\bullet}=A_{\sim p,\bullet}.
\]
This proves Proposition \ref{prop.addcol.props1} \textbf{(f)}.
\end{proof}

\begin{proof}
[Proof of Proposition \ref{prop.addcol.props2}.]We have $n-1\in\mathbb{N}$
(since $n$ is a positive integer).

\textbf{(a)} Write the $n\times\left(  n-1\right)  $-matrix $A$ in the form
$A=\left(  a_{i,j}\right)  _{1\leq i\leq n,\ 1\leq j\leq n-1}\in
\mathbb{K}^{n\times\left(  n-1\right)  }$. Let $v=\left(  v_{1},v_{2}%
,\ldots,v_{n}\right)  ^{T}\in\mathbb{K}^{n\times1}$ be a vector.

Clearly, $\left(  A\mid v\right)  $ is an $n\times\left(  \left(  n-1\right)
+1\right)  $-matrix (since $A$ is an $n\times\left(  n-1\right)  $-matrix, and
since $v$ is an $n\times1$-matrix). In other words, $\left(  A\mid v\right)  $
is an $n\times n$-matrix (since $\left(  n-1\right)  +1=n$). Thus,
$\det\left(  A\mid v\right)  $ is well-defined.

Write the matrix $\left(  A\mid v\right)  $ in the form $\left(  A\mid
v\right)  =\left(  b_{i,j}\right)  _{1\leq i\leq n,\ 1\leq j\leq n}$. (This is
possible since $\left(  A\mid v\right)  $ is an $n\times n$-matrix.) Thus,
\[
\left(  A\mid v\right)  =\left(  b_{i,j}\right)  _{1\leq i\leq n,\ 1\leq j\leq
n}=\left(  b_{i,j}\right)  _{1\leq i\leq n,\ 1\leq j\leq\left(  n-1\right)
+1}%
\]
(since $n=\left(  n-1\right)  +1$). Hence,
(\ref{eq.lem.sol.prop.addcol.props.1.2}) shows that we have%
\[
b_{i,\left(  n-1\right)  +1}=v_{i}%
\]
for every $i\in\left\{  1,2,\ldots,n\right\}  $. Since $\left(  n-1\right)
+1=n$, this rewrites as follows: We have%
\begin{equation}
b_{i,n}=v_{i} \label{pf.prop.addcols.props2.1}%
\end{equation}
for every $i\in\left\{  1,2,\ldots,n\right\}  $.

But every $p\in\left\{  1,2,\ldots,n\right\}  $ satisfies%
\begin{equation}
\left(  A\mid v\right)  _{\sim p,\sim n}=A_{\sim p,\bullet}
\label{pf.prop.addcols.props2.2}%
\end{equation}
\footnote{\textit{Proof of (\ref{pf.prop.addcols.props2.2}):} Let
$p\in\left\{  1,2,\ldots,n\right\}  $. Proposition \ref{prop.addcol.props1}
\textbf{(f)} (applied to $m=n-1$) yields $\left(  A\mid v\right)  _{\sim
p,\sim\left(  \left(  n-1\right)  +1\right)  }=A_{\sim p,\bullet}$. This
rewrites as $\left(  A\mid v\right)  _{\sim p,\sim n}=A_{\sim p,\bullet}$
(since $\left(  n-1\right)  +1=n$). This proves
(\ref{pf.prop.addcols.props2.2}).}.

We have $n\in\left\{  1,2,\ldots,n\right\}  $ (since $n$ is a positive
integer) and $\left(  A\mid v\right)  =\left(  b_{i,j}\right)  _{1\leq i\leq
n,\ 1\leq j\leq n}$. Hence, Theorem \ref{thm.laplace.gen} \textbf{(b)}
(applied to $\left(  A\mid v\right)  $, $b_{i,j}$ and $n$ instead of $A$,
$a_{i,j}$ and $q$) yields%
\begin{align*}
\det\left(  A\mid v\right)   &  =\sum_{p=1}^{n}\underbrace{\left(  -1\right)
^{p+n}}_{\substack{=\left(  -1\right)  ^{n+p}\\\text{(since }p+n=n+p\text{)}%
}}\underbrace{b_{p,n}}_{\substack{=v_{p}\\\text{(by
(\ref{pf.prop.addcols.props2.1}),}\\\text{applied to }i=p\text{)}}}\det\left(
\underbrace{\left(  A\mid v\right)  _{\sim p,\sim n}}_{\substack{=A_{\sim
p,\bullet}\\\text{(by (\ref{pf.prop.addcols.props2.2}))}}}\right) \\
&  =\sum_{p=1}^{n}\left(  -1\right)  ^{n+p}v_{p}\det\left(  A_{\sim p,\bullet
}\right)  =\sum_{i=1}^{n}\left(  -1\right)  ^{n+i}v_{i}\det\left(  A_{\sim
i,\bullet}\right)
\end{align*}
(here, we have renamed the summation index $p$ as $i$). This proves
Proposition \ref{prop.addcol.props2} \textbf{(a)}.

\textbf{(b)} Let $p\in\left\{  1,2,\ldots,n\right\}  $. For any two objects
$i$ and $j$, we define an element $\delta_{i,j}\in\mathbb{K}$ by $\delta
_{i,j}=%
\begin{cases}
1, & \text{if }i=j;\\
0, & \text{if }i\neq j
\end{cases}
$. Then, $I_{n}=\left(  \delta_{i,j}\right)  _{1\leq i\leq n,\ 1\leq j\leq n}$
(by the definition of $I_{n}$). Now, $\left(  I_{n}\right)  _{\bullet,p}$ is
the $p$-th column of the matrix $I_{n}$ (by the definition of $\left(
I_{n}\right)  _{\bullet,p}$). Thus,%
\begin{align*}
\left(  I_{n}\right)  _{\bullet,p}  &  =\left(  \text{the }p\text{-th column
of the matrix }I_{n}\right) \\
&  =\left(
\begin{array}
[c]{c}%
\delta_{1,p}\\
\delta_{2,p}\\
\vdots\\
\delta_{n,p}%
\end{array}
\right)  \ \ \ \ \ \ \ \ \ \ \left(  \text{since }I_{n}=\left(  \delta
_{i,j}\right)  _{1\leq i\leq n,\ 1\leq j\leq n}\right) \\
&  =\left(  \delta_{1,p},\delta_{2,p},\ldots,\delta_{n,p}\right)  ^{T}.
\end{align*}
Hence, Proposition \ref{prop.addcol.props2} \textbf{(a)} (applied to $\left(
I_{n}\right)  _{\bullet,p}$ and $\delta_{i,p}$ instead of $v$ and $v_{i}$)
yields%
\begin{align*}
&  \det\left(  A\mid\left(  I_{n}\right)  _{\bullet,p}\right) \\
&  =\underbrace{\sum_{i=1}^{n}}_{=\sum_{i\in\left\{  1,2,\ldots,n\right\}  }%
}\left(  -1\right)  ^{n+i}\delta_{i,p}\det\left(  A_{\sim i,\bullet}\right)
=\sum_{i\in\left\{  1,2,\ldots,n\right\}  }\left(  -1\right)  ^{n+i}%
\delta_{i,p}\det\left(  A_{\sim i,\bullet}\right) \\
&  =\left(  -1\right)  ^{n+p}\underbrace{\delta_{p,p}}%
_{\substack{=1\\\text{(since }p=p\text{)}}}\det\left(  A_{\sim p,\bullet
}\right)  +\sum_{\substack{i\in\left\{  1,2,\ldots,n\right\}  ;\\i\neq
p}}\left(  -1\right)  ^{n+i}\underbrace{\delta_{i,p}}%
_{\substack{=0\\\text{(since }i\neq p\text{)}}}\det\left(  A_{\sim i,\bullet
}\right) \\
&  \ \ \ \ \ \ \ \ \ \ \left(
\begin{array}
[c]{c}%
\text{here, we have split off the addend for }i=p\text{ from the sum,}\\
\text{since }p\in\left\{  1,2,\ldots,n\right\}
\end{array}
\right) \\
&  =\left(  -1\right)  ^{n+p}\det\left(  A_{\sim p,\bullet}\right)
+\underbrace{\sum_{\substack{i\in\left\{  1,2,\ldots,n\right\}  ;\\i\neq
p}}\left(  -1\right)  ^{n+i}0\det\left(  A_{\sim i,\bullet}\right)  }%
_{=0}=\left(  -1\right)  ^{n+p}\det\left(  A_{\sim p,\bullet}\right)  .
\end{align*}
This proves Proposition \ref{prop.addcol.props2} \textbf{(b)}.
\end{proof}

Before we prove Proposition \ref{prop.addcol.props3}, let us make two more
trivial observations:

\begin{lemma}
\label{lem.sol.prop.addcol.props.cols}Let $n\in\mathbb{N}$ and $m\in
\mathbb{N}$. Let $A$ and $B$ be two $n\times m$-matrices. Assume that
$A_{\bullet,q}=B_{\bullet,q}$ for each $q\in\left\{  1,2,\ldots,m\right\}  $.
Then, $A=B$.
\end{lemma}

\begin{proof}
[Proof of Lemma \ref{lem.sol.prop.addcol.props.cols}.]We have assumed that%
\begin{equation}
A_{\bullet,q}=B_{\bullet,q} \label{pf.lem.sol.prop.addcol.props.cols.AB}%
\end{equation}
for each $q\in\left\{  1,2,\ldots,m\right\}  $. For every $q\in\left\{
1,2,\ldots,m\right\}  $, we have%
\begin{equation}
A_{\bullet,q}=\left(  \text{the }q\text{-th column of the matrix }A\right)
\label{pf.lem.sol.prop.addcol.props.cols.A}%
\end{equation}
(since $A_{\bullet,q}$ is the $q$-th column of the matrix $A$ (by the
definition of $A_{\bullet,q}$)) and%
\begin{equation}
B_{\bullet,q}=\left(  \text{the }q\text{-th column of the matrix }B\right)
\label{pf.lem.sol.prop.addcol.props.cols.B}%
\end{equation}
(since $B_{\bullet,q}$ is the $q$-th column of the matrix $B$ (by the
definition of $B_{\bullet,q}$)). Hence, for every $q\in\left\{  1,2,\ldots
,m\right\}  $, we have%
\begin{align}
&  \left(  \text{the }q\text{-th column of the matrix }A\right) \nonumber\\
&  =A_{\bullet,q}\ \ \ \ \ \ \ \ \ \ \left(  \text{by
(\ref{pf.lem.sol.prop.addcol.props.cols.A})}\right) \nonumber\\
&  =B_{\bullet,q}\ \ \ \ \ \ \ \ \ \ \left(  \text{by
(\ref{pf.lem.sol.prop.addcol.props.cols.AB})}\right) \nonumber\\
&  =\left(  \text{the }q\text{-th column of the matrix }B\right)
\ \ \ \ \ \ \ \ \ \ \left(  \text{by
(\ref{pf.lem.sol.prop.addcol.props.cols.B})}\right)  .
\label{pf.lem.sol.prop.addcol.props.cols.q=q}%
\end{align}

Write the matrix $A$ in the form $A=\left(  a_{i,j}\right)  _{1\leq i\leq
n,\ 1\leq j\leq m}$. (This is possible since $A$ is an $n\times m$-matrix.)

Write the matrix $B$ in the form $B=\left(  b_{i,j}\right)  _{1\leq i\leq
n,\ 1\leq j\leq m}$. (This is possible since $B$ is an $n\times m$-matrix.)

We have $A=\left(  a_{i,j}\right)  _{1\leq i\leq n,\ 1\leq j\leq m}$. Thus,%
\begin{equation}
\left(  \text{the }\left(  i,j\right)  \text{-th entry of the matrix
}A\right)  =a_{i,j} \label{pf.lem.sol.prop.addcol.props.cols.Aij}%
\end{equation}
for every $i\in\left\{  1,2,\ldots,n\right\}  $ and $j\in\left\{
1,2,\ldots,m\right\}  $.

We have $B=\left(  b_{i,j}\right)  _{1\leq i\leq n,\ 1\leq j\leq m}$. Thus,%
\begin{equation}
\left(  \text{the }\left(  i,j\right)  \text{-th entry of the matrix
}B\right)  =b_{i,j} \label{pf.lem.sol.prop.addcol.props.cols.Bij}%
\end{equation}
for every $i\in\left\{  1,2,\ldots,n\right\}  $ and $j\in\left\{
1,2,\ldots,m\right\}  $.

Now, for every $i\in\left\{  1,2,\ldots,n\right\}  $ and $j\in\left\{
1,2,\ldots,m\right\}  $, we have%
\begin{align*}
a_{i,j}  &  =\left(  \text{the }\left(  i,j\right)  \text{-th entry of the
matrix }A\right)  \ \ \ \ \ \ \ \ \ \ \left(  \text{by
(\ref{pf.lem.sol.prop.addcol.props.cols.Aij})}\right) \\
&  =\left(  \text{the }i\text{-th entry of }\underbrace{\text{the }j\text{-th
column of the matrix }A}_{\substack{=\left(  \text{the }j\text{-th column of
the matrix }B\right)  \\\text{(by (\ref{pf.lem.sol.prop.addcol.props.cols.q=q}%
), applied to }q=j\text{)}}}\right) \\
&  =\left(  \text{the }i\text{-th entry of the }j\text{-th column of the
matrix }B\right) \\
&  =\left(  \text{the }\left(  i,j\right)  \text{-th entry of the matrix
}B\right) \\
&  =b_{i,j}\ \ \ \ \ \ \ \ \ \ \left(  \text{by
(\ref{pf.lem.sol.prop.addcol.props.cols.Bij})}\right)  .
\end{align*}
Hence, $\left(  \underbrace{a_{i,j}}_{=b_{i,j}}\right)  _{1\leq i\leq
n,\ 1\leq j\leq m}=\left(  b_{i,j}\right)  _{1\leq i\leq n,\ 1\leq j\leq m}=B$
(since $B=\left(  b_{i,j}\right)  _{1\leq i\leq n,\ 1\leq j\leq m}$). Thus,
$A=\left(  a_{i,j}\right)  _{1\leq i\leq n,\ 1\leq j\leq m}=B$. This proves
Lemma \ref{lem.sol.prop.addcol.props.cols}.
\end{proof}

\begin{lemma}
\label{lem.sol.prop.addcol.props3.bc}Let $n\in\mathbb{N}$. Let $A=\left(
a_{i,j}\right)  _{1\leq i\leq n,\ 1\leq j\leq n}\in\mathbb{K}^{n\times n}$ be
an $n\times n$-matrix. Let $r$ and $q$ be two elements of $\left\{
1,2,\ldots,n\right\}  $. Then,%
\begin{equation}
\det\left(  A_{\bullet,\sim q}\mid A_{\bullet,r}\right)  =\left(  -1\right)
^{n+q}\sum_{p=1}^{n}\left(  -1\right)  ^{p+q}a_{p,r}\det\left(  A_{\sim p,\sim
q}\right)  .\nonumber
\end{equation}

\end{lemma}

\begin{proof}
[Proof of Lemma \ref{lem.sol.prop.addcol.props3.bc}.]We have $r\in\left\{
1,2,\ldots,n\right\}  $. Hence, $A_{\bullet,r}$ is an $n\times1$-matrix (since
$A$ is an $n\times n$-matrix). Also, $r\in\left\{  1,2,\ldots,n\right\}  $
shows that $1\leq r\leq n$; hence, $1\leq n$, so that $n$ is a positive integer.

Also, $q\in\left\{  1,2,\ldots,n\right\}  $. Hence, $A_{\bullet,\sim q}$ is an
$n\times\left(  n-1\right)  $-matrix (since $A$ is an $n\times n$-matrix); in
other words, $A_{\bullet,\sim q}\in\mathbb{K}^{n\times\left(  n-1\right)  }$.

Furthermore, $A_{\bullet,r}$ is the $r$-th column of the matrix $A$ (by the
definition of $A_{\bullet,r}$). Thus,%
\begin{align*}
A_{\bullet,r}  &  =\left(  \text{the }r\text{-th column of the matrix
}A\right) \\
&  =\left(
\begin{array}
[c]{c}%
a_{1,r}\\
a_{2,r}\\
\vdots\\
a_{n,r}%
\end{array}
\right)  \ \ \ \ \ \ \ \ \ \ \left(  \text{since }A=\left(  a_{i,j}\right)
_{1\leq i\leq n,\ 1\leq j\leq n}\right) \\
&  =\left(  a_{1,r},a_{2,r},\ldots,a_{n,r}\right)  ^{T}.
\end{align*}
Thus, Proposition \ref{prop.addcol.props2} \textbf{(a)} (applied to
$A_{\bullet,\sim q}$, $A_{\bullet,r}$ and $a_{i,r}$ instead of $A$, $v$ and
$v_{i}$) yields%
\begin{equation}
\det\left(  A_{\bullet,\sim q}\mid A_{\bullet,r}\right)  =\sum_{i=1}%
^{n}\left(  -1\right)  ^{n+i}a_{i,r}\det\left(  \left(  A_{\bullet,\sim
q}\right)  _{\sim i,\bullet}\right)  .
\label{pf.lem.sol.prop.addcol.props3.bc.1}%
\end{equation}
But every $i\in\left\{  1,2,\ldots,n\right\}  $ satisfies%
\begin{equation}
\left(  A_{\bullet,\sim q}\right)  _{\sim i,\bullet}=A_{\sim i,\sim q}
\label{pf.lem.sol.prop.addcol.props3.bc.2}%
\end{equation}
\footnote{\textit{Proof of (\ref{pf.lem.sol.prop.addcol.props3.bc.2}):} Let
$i\in\left\{  1,2,\ldots,n\right\}  $. Then, Proposition
\ref{prop.unrows.basics} \textbf{(c)} (applied to $n$, $i$ and $q$ instead of
$m$, $u$ and $v$) yields $\left(  A_{\bullet,\sim q}\right)  _{\sim i,\bullet
}=\left(  A_{\sim i,\bullet}\right)  _{\bullet,\sim q}=A_{\sim i,\sim q}$.
Hence, $\left(  A_{\bullet,\sim q}\right)  _{\sim i,\bullet}=A_{\sim i,\sim
q}$. This proves (\ref{pf.lem.sol.prop.addcol.props3.bc.2}).} and%
\begin{equation}
\left(  -1\right)  ^{n+i}=\left(  -1\right)  ^{n+q}\left(  -1\right)  ^{i+q}
\label{pf.lem.sol.prop.addcol.props3.bc.3}%
\end{equation}
\footnote{\textit{Proof of (\ref{pf.lem.sol.prop.addcol.props3.bc.3}):} Let
$i\in\left\{  1,2,\ldots,n\right\}  $. Then, $\left(  -1\right)  ^{n+q}\left(
-1\right)  ^{i+q}=\left(  -1\right)  ^{\left(  n+q\right)  +\left(
i+q\right)  }=\left(  -1\right)  ^{n+i}$ (since $\left(  n+q\right)  +\left(
i+q\right)  =n+i+2q\equiv n+i\operatorname{mod}2$). This proves
(\ref{pf.lem.sol.prop.addcol.props3.bc.3}).}. Now,
(\ref{pf.lem.sol.prop.addcol.props3.bc.1}) becomes%
\begin{align*}
\det\left(  A_{\bullet,\sim q}\mid A_{\bullet,r}\right)   &  =\sum_{i=1}%
^{n}\underbrace{\left(  -1\right)  ^{n+i}}_{\substack{=\left(  -1\right)
^{n+q}\left(  -1\right)  ^{i+q}\\\text{(by
(\ref{pf.lem.sol.prop.addcol.props3.bc.3}))}}}a_{i,r}\det\left(
\underbrace{\left(  A_{\bullet,\sim q}\right)  _{\sim i,\bullet}%
}_{\substack{=A_{\sim i,\sim q}\\\text{(by
(\ref{pf.lem.sol.prop.addcol.props3.bc.2}))}}}\right) \\
&  =\sum_{i=1}^{n}\left(  -1\right)  ^{n+q}\left(  -1\right)  ^{i+q}%
a_{i,r}\det\left(  A_{\sim i,\sim q}\right) \\
&  =\left(  -1\right)  ^{n+q}\sum_{i=1}^{n}\left(  -1\right)  ^{i+q}%
a_{i,r}\det\left(  A_{\sim i,\sim q}\right) \\
&  =\left(  -1\right)  ^{n+q}\sum_{p=1}^{n}\left(  -1\right)  ^{p+q}%
a_{p,r}\det\left(  A_{\sim p,\sim q}\right) \\
&  \ \ \ \ \ \ \ \ \ \ \left(  \text{here, we have renamed the summation index
}i\text{ as }p\right)  .
\end{align*}
This proves Lemma \ref{lem.sol.prop.addcol.props3.bc}.
\end{proof}

\begin{proof}
[Proof of Proposition \ref{prop.addcol.props3}.]\textbf{(a)} Assume that
$n>0$. Thus, $n$ is a positive integer, so that $n\in\left\{  1,2,\ldots
,n\right\}  $. Hence, $A_{\bullet,\sim n}$ is a well-defined $n\times\left(
n-1\right)  $-matrix, and $A_{\bullet,n}$ is a well-defined $n\times1$-matrix.
Therefore, $\left(  A_{\bullet,\sim n}\mid A_{\bullet,n}\right)  $ is an
$n\times\left(  \left(  n-1\right)  +1\right)  $-matrix. In other words,
$\left(  A_{\bullet,\sim n}\mid A_{\bullet,n}\right)  $ is an $n\times
n$-matrix (since $\left(  n-1\right)  +1=n$). Also, $n-1\in\mathbb{N}$ (since
$n$ is a positive integer).

We shall now show that%
\begin{equation}
A_{\bullet,q}=\left(  A_{\bullet,\sim n}\mid A_{\bullet,n}\right)
_{\bullet,q} \label{pf.prop.addcol.props3.a.main}%
\end{equation}
for each $q\in\left\{  1,2,\ldots,n\right\}  $.

[\textit{Proof of (\ref{pf.prop.addcol.props3.a.main}):} Let $q\in\left\{
1,2,\ldots,n\right\}  $. We must prove the equality
(\ref{pf.prop.addcol.props3.a.main}). We are in one of the following two cases:

\textit{Case 1:} We have $q\neq n$.

\textit{Case 2:} We have $q=n$.

Let us first consider Case 1. In this case, we have $q\neq n$. Combining
$q\in\left\{  1,2,\ldots,n\right\}  $ with $q\neq n$, we obtain $q\in\left\{
1,2,\ldots,n\right\}  \setminus\left\{  n\right\}  =\left\{  1,2,\ldots
,n-1\right\}  $. Thus, Proposition \ref{prop.addcol.props1} \textbf{(a)}
(applied to $n-1$, $A_{\bullet,\sim n}$ and $A_{\bullet,n}$ instead of $m$,
$A$ and $v$) yields $\left(  A_{\bullet,\sim n}\mid A_{\bullet,n}\right)
_{\bullet,q}=\left(  A_{\bullet,\sim n}\right)  _{\bullet,q}$.

But Proposition \ref{prop.unrows.basics} \textbf{(d)} (applied to $n$, $n$ and
$q$ instead of $m$, $v$ and $w$) yields $\left(  A_{\bullet,\sim n}\right)
_{\bullet,q}=A_{\bullet,q}$. Hence, $A_{\bullet,q}=\left(  A_{\bullet,\sim
n}\right)  _{\bullet,q}=\left(  A_{\bullet,\sim n}\mid A_{\bullet,n}\right)
_{\bullet,q}$ (since $\left(  A_{\bullet,\sim n}\mid A_{\bullet,n}\right)
_{\bullet,q}=\left(  A_{\bullet,\sim n}\right)  _{\bullet,q}$). Hence,
(\ref{pf.prop.addcol.props3.a.main}) is proven in Case 1.

Let us now consider Case 2. In this case, we have $q=n$. But Proposition
\ref{prop.addcol.props1} \textbf{(b)} (applied to $n-1$, $A_{\bullet,\sim n}$
and $A_{\bullet,n}$ instead of $m$, $A$ and $v$) yields $\left(
A_{\bullet,\sim n}\mid A_{\bullet,n}\right)  _{\bullet,\left(  n-1\right)
+1}=A_{\bullet,n}$. Thus, $A_{\bullet,n}=\left(  A_{\bullet,\sim n}\mid
A_{\bullet,n}\right)  _{\bullet,\left(  n-1\right)  +1}=\left(  A_{\bullet
,\sim n}\mid A_{\bullet,n}\right)  _{\bullet,q}$ (since $\left(  n-1\right)
+1=n=q$). Now, $q=n$, so that $A_{\bullet,q}=A_{\bullet,n}=\left(
A_{\bullet,\sim n}\mid A_{\bullet,n}\right)  _{\bullet,q}$. Hence,
(\ref{pf.prop.addcol.props3.a.main}) is proven in Case 2.

Now, (\ref{pf.prop.addcol.props3.a.main}) is proven in each of the two Cases 1
and 2. Since these two Cases cover all possibilities, this completes the proof
of (\ref{pf.prop.addcol.props3.a.main}).]

Now we know that $A$ and $\left(  A_{\bullet,\sim n}\mid A_{\bullet,n}\right)
$ are two $n\times n$-matrices and that (\ref{pf.prop.addcol.props3.a.main})
holds for each $q\in\left\{  1,2,\ldots,n\right\}  $. Hence, Lemma
\ref{lem.sol.prop.addcol.props.cols} (applied to $m=n$ and $B=\left(
A_{\bullet,\sim n}\mid A_{\bullet,n}\right)  $) yields $A=\left(
A_{\bullet,\sim n}\mid A_{\bullet,n}\right)  $. This proves Proposition
\ref{prop.addcol.props3} \textbf{(a)}.

\textbf{(b)} Write the matrix $A$ in the form $A=\left(  a_{i,j}\right)
_{1\leq i\leq n,\ 1\leq j\leq n}$. (This is possible since $A$ is an $n\times
n$-matrix.)

Lemma \ref{lem.sol.prop.addcol.props3.bc} (applied to $r=q$) yields%
\begin{equation}
\det\left(  A_{\bullet,\sim q}\mid A_{\bullet,q}\right)  =\left(  -1\right)
^{n+q}\sum_{p=1}^{n}\left(  -1\right)  ^{p+q}a_{p,q}\det\left(  A_{\sim p,\sim
q}\right)  .\nonumber
\end{equation}
Comparing this with%
\[
\left(  -1\right)  ^{n+q}\underbrace{\det A}_{\substack{=\sum_{p=1}^{n}\left(
-1\right)  ^{p+q}a_{p,q}\det\left(  A_{\sim p,\sim q}\right)  \\\text{(by
Theorem \ref{thm.laplace.gen} \textbf{(b)})}}}=\left(  -1\right)  ^{n+q}%
\sum_{p=1}^{n}\left(  -1\right)  ^{p+q}a_{p,q}\det\left(  A_{\sim p,\sim
q}\right)  ,
\]
we obtain $\det\left(  A_{\bullet,\sim q}\mid A_{\bullet,q}\right)  =\left(
-1\right)  ^{n+q}\det A$. This proves Proposition \ref{prop.addcol.props3}
\textbf{(b)}.

[\textit{Remark:} Let me outline an alternative proof of Proposition
\ref{prop.addcol.props3} \textbf{(b)}: Let $q\in\left\{  1,2,\ldots,n\right\}
$. The matrix $\left(  A_{\bullet,\sim q}\mid A_{\bullet,q}\right)  $ is
obtained from the matrix $A$ by removing the $q$-th column and then
reattaching this column to the right end of the matrix. This procedure can be
replaced by the following procedure, which clearly leads to the same result:

\begin{itemize}
\item Swap the $q$-th column of $A$ with the $\left(  q+1\right)  $-th column;

\item then swap the $\left(  q+1\right)  $-th column of the resulting matrix
with the $\left(  q+2\right)  $-th column;

\item then swap the $\left(  q+2\right)  $-th column of the resulting matrix
with the $\left(  q+3\right)  $-th column;

\item and so on, finally swapping the $\left(  n-1\right)  $-st column of the
matrix with the $n$-th column.
\end{itemize}

But this latter procedure is a sequence of $n-q$ swaps of two columns. Each
such swap multiplies the determinant of the matrix by $-1$ (according to
Exercise \ref{exe.ps4.6} \textbf{(b)}). Thus, the whole procedure multiplies
the determinant of the matrix by $\left(  -1\right)  ^{n-q}=\left(  -1\right)
^{n+q}$ (since $n-q\equiv n+q\operatorname{mod}2$). Since this procedure takes
the matrix $A$ to the matrix $\left(  A_{\bullet,\sim q}\mid A_{\bullet
,q}\right)  $, we thus conclude that $\det\left(  A_{\bullet,\sim q}\mid
A_{\bullet,q}\right)  =\left(  -1\right)  ^{n+q}\det A$. This proves
Proposition \ref{prop.addcol.props3} \textbf{(b)} again.]

\textbf{(c)} Let $r$ and $q$ be two elements of $\left\{  1,2,\ldots
,n\right\}  $ satisfying $r\neq q$. We must show that $\det\left(
A_{\bullet,\sim q}\mid A_{\bullet,r}\right)  =0$.

Write the matrix $A$ in the form $A=\left(  a_{i,j}\right)  _{1\leq i\leq
n,\ 1\leq j\leq n}$. (This is possible since $A$ is an $n\times n$-matrix.)

We have $q\neq r$ (since $r\neq q$). Hence, Proposition \ref{prop.laplace.0}
\textbf{(b)} yields%
\begin{equation}
0=\sum_{p=1}^{n}\left(  -1\right)  ^{p+q}a_{p,r}\det\left(  A_{\sim p,\sim
q}\right)  . \label{pf.prop.addcol.props3.c.1}%
\end{equation}

Lemma \ref{lem.sol.prop.addcol.props3.bc} yields%
\begin{equation}
\det\left(  A_{\bullet,\sim q}\mid A_{\bullet,r}\right)  =\left(  -1\right)
^{n+q}\underbrace{\sum_{p=1}^{n}\left(  -1\right)  ^{p+q}a_{p,r}\det\left(
A_{\sim p,\sim q}\right)  }_{\substack{=0\\\text{(by
(\ref{pf.prop.addcol.props3.c.1}))}}}=0.\nonumber
\end{equation}
This proves Proposition \ref{prop.addcol.props3} \textbf{(c)}.

[\textit{Remark:} Let me outline an alternative proof of Proposition
\ref{prop.addcol.props3} \textbf{(c)}: Let $r$ and $q$ be two elements of
$\left\{  1,2,\ldots,n\right\}  $ satisfying $r\neq q$. Then, $A_{\bullet,r}$
is one of the columns of the matrix $A_{\bullet,\sim q}$ (since the $r$-th
column of the matrix $A$ is not lost when the $q$-th column is removed).
Hence, the column vector $A_{\bullet,r}$ appears twice as a column in the
matrix $\left(  A_{\bullet,\sim q}\mid A_{\bullet,r}\right)  $ (namely, once
as one of the columns of $A_{\bullet,\sim q}$, and another time as the
attached column). Therefore, the matrix $\left(  A_{\bullet,\sim q}\mid
A_{\bullet,r}\right)  $ has two equal columns. Exercise \ref{exe.ps4.6}
\textbf{(f)} thus shows that $\det\left(  A_{\bullet,\sim q}\mid A_{\bullet
,r}\right)  =0$. This proves Proposition \ref{prop.addcol.props3} \textbf{(c)} again.]

\textbf{(d)} Let $p\in\left\{  1,2,\ldots,n\right\}  $ and $q\in\left\{
1,2,\ldots,n\right\}  $. Then, Proposition \ref{prop.unrows.basics}
\textbf{(c)} (applied to $u=p$ and $v=q$) yields%
\[
\left(  A_{\bullet,\sim q}\right)  _{\sim p,\bullet}=\left(  A_{\sim
p,\bullet}\right)  _{\bullet,\sim q}=A_{\sim p,\sim q}.
\]

But $A$ is an $n\times n$-matrix. Hence, $A_{\bullet,\sim q}$ is an
$n\times\left(  n-1\right)  $-matrix (since $q\in\left\{  1,2,\ldots
,n\right\}  $). Moreover, $p\in\left\{  1,2,\ldots,n\right\}  $, so that
$1\leq p\leq n$ and thus $n\geq1$; hence, $n$ is a positive integer.
Proposition \ref{prop.addcol.props2} \textbf{(b)} (applied to $A_{\bullet,\sim
q}$ instead of $A$) thus yields
\[
\det\left(  A_{\bullet,\sim q}\mid\left(  I_{n}\right)  _{\bullet,p}\right)
=\left(  -1\right)  ^{n+p}\det\left(  \underbrace{\left(  A_{\bullet,\sim
q}\right)  _{\sim p,\bullet}}_{=A_{\sim p,\sim q}}\right)  =\left(  -1\right)
^{n+p}\det\left(  A_{\sim p,\sim q}\right)  .
\]
This proves Proposition \ref{prop.addcol.props3} \textbf{(d)}.

\textbf{(e)} Let $u$ and $v$ be two elements of $\left\{  1,2,\ldots
,n\right\}  $ satisfying $u<v$. Let $r$ be an element of $\left\{
1,2,\ldots,n-1\right\}  $ satisfying $r\neq u$. We must show that $\det\left(
A_{\bullet,\sim u}\mid\left(  A_{\bullet,\sim v}\right)  _{\bullet,r}\right)
=0$.

We are in one of the following two cases:

\textit{Case 1:} We have $r<v$.

\textit{Case 2:} We have $r\geq v$.

Let us first consider Case 1. In this case, we have $r<v$. Thus, $r\leq v-1$
(since $r$ and $v$ are integers). But $r\geq1$ (since $r\in\left\{
1,2,\ldots,n-1\right\}  $). Combining $r\geq1$ with $r\leq v-1$, we obtain
$r\in\left\{  1,2,\ldots,v-1\right\}  $. Hence, Proposition
\ref{prop.unrows.basics} \textbf{(d)} (applied to $m=n$ and $w=r$) yields
$\left(  A_{\bullet,\sim v}\right)  _{\bullet,r}=A_{\bullet,r}$. Hence,%
\[
\det\left(  A_{\bullet,\sim u}\mid\underbrace{\left(  A_{\bullet,\sim
v}\right)  _{\bullet,r}}_{=A_{\bullet,r}}\right)  =\det\left(  A_{\bullet,\sim
u}\mid A_{\bullet,r}\right)  =0
\]
(by Proposition \ref{prop.addcol.props3} \textbf{(c)}, applied to $q=u$).
Hence, Proposition \ref{prop.addcol.props3} \textbf{(e)} is proven in Case 1.

Let us now consider Case 2. In this case, we have $r\geq v$. But $r\leq n-1$
(since $r\in\left\{  1,2,\ldots,n-1\right\}  $), so that $r+1\leq n$.
Combining this with $r+1\geq r\geq1$ (since $r\in\left\{  1,2,\ldots
,n\right\}  $), we obtain $r+1\in\left\{  1,2,\ldots,n\right\}  $.
Furthermore, $r+1>r\geq v>u$ (since $u<v$) and thus $r+1\neq u$. Hence,
Proposition \ref{prop.addcol.props3} \textbf{(c)} (applied to $r+1$ and $u$
instead of $r$ and $q$) yields $\det\left(  A_{\bullet,\sim u}\mid
A_{\bullet,r+1}\right)  =0$.

But combining $r\geq v$ with $r\leq n-1$, we obtain $r\in\left\{
v,v+1,\ldots,n-1\right\}  $. Thus, Proposition \ref{prop.unrows.basics}
\textbf{(e)} (applied to $m=n$ and $w=r$) yields $\left(  A_{\bullet,\sim
v}\right)  _{\bullet,r}=A_{\bullet,r+1}$. Hence,%
\[
\det\left(  A_{\bullet,\sim u}\mid\underbrace{\left(  A_{\bullet,\sim
v}\right)  _{\bullet,r}}_{=A_{\bullet,r+1}}\right)  =\det\left(
A_{\bullet,\sim u}\mid A_{\bullet,r+1}\right)  =0.
\]
Hence, Proposition \ref{prop.addcol.props3} \textbf{(e)} is proven in Case 2.

We have now proven Proposition \ref{prop.addcol.props3} \textbf{(e)} in each
of the two Cases 1 and 2. Since these two Cases cover all possibilities, this
shows that Proposition \ref{prop.addcol.props3} \textbf{(e)} always holds.

\textbf{(f)} Let $u$ and $v$ be two elements of $\left\{  1,2,\ldots
,n\right\}  $ satisfying $u<v$. Then, $u<v$, so that $u\leq v-1$ (since $u$
and $v$ are integers). Combined with $u\geq1$ (since $u\in\left\{
1,2,\ldots,n\right\}  $), this yields $u\in\left\{  1,2,\ldots,v-1\right\}  $.
Hence, Proposition \ref{prop.unrows.basics} \textbf{(d)} (applied to $m=n$ and
$w=u$) yields $\left(  A_{\bullet,\sim v}\right)  _{\bullet,u}=A_{\bullet,u}$.
Thus,%
\[
\det\left(  A_{\bullet,\sim u}\mid\underbrace{\left(  A_{\bullet,\sim
v}\right)  _{\bullet,u}}_{=A_{\bullet,u}}\right)  =\det\left(  A_{\bullet,\sim
u}\mid A_{\bullet,u}\right)  =\left(  -1\right)  ^{n+u}\det A
\]
(by Proposition \ref{prop.addcol.props3} \textbf{(b)}, applied to $q=u$).
Thus,%
\[
\left(  -1\right)  ^{u}\underbrace{\det\left(  A_{\bullet,\sim u}\mid\left(
A_{\bullet,\sim v}\right)  _{\bullet,u}\right)  }_{=\left(  -1\right)
^{n+u}\det A}=\underbrace{\left(  -1\right)  ^{u}\left(  -1\right)  ^{n+u}%
}_{\substack{=\left(  -1\right)  ^{u+\left(  n+u\right)  }=\left(  -1\right)
^{n}\\\text{(since }u+\left(  n+u\right)  =2u+n\equiv n\operatorname{mod}%
2\text{)}}}\det A=\left(  -1\right)  ^{n}\det A.
\]
This proves Proposition \ref{prop.addcol.props3} \textbf{(f)}.
\end{proof}
\end{verlong}

\begin{proof}
[Solution to Exercise \ref{exe.prop.addcol.props}.]We have proven Proposition
\ref{prop.addcol.props1}, Proposition \ref{prop.addcol.props2} and Proposition
\ref{prop.addcol.props3}. Thus, Exercise \ref{exe.prop.addcol.props} is solved.
\end{proof}

\subsection{Solution to Exercise \ref{exe.desnanot.jaw}}

\begin{proof}
[Solution to Exercise \ref{exe.desnanot.jaw}.]We have $u<v$ and $u<w$. Define
an $n\times n$-matrix $C\in\mathbb{K}^{n\times n}$ as in Lemma
\ref{lem.desnanot.AB.tech}.

Lemma \ref{lem.desnanot.AB.tech} \textbf{(a)} yields%
\begin{align}
\det\left(  C_{\sim v,\sim q}\right)   &  =-\left(  -1\right)  ^{n+u}%
\underbrace{\det\left(  \operatorname*{rows}\nolimits_{1,2,\ldots
,\widehat{u},\ldots,\widehat{v},\ldots,n}\left(  B_{\bullet,\sim q}\right)
\right)  }_{\substack{=\beta_{u,v}\\\text{(since }\beta_{u,v}=\det\left(
\operatorname*{rows}\nolimits_{1,2,\ldots,\widehat{u},\ldots,\widehat{v}%
,\ldots,n}\left(  B_{\bullet,\sim q}\right)  \right)  \\\text{(by the
definition of }\beta_{u,v}\text{))}}}\nonumber\\
&  =-\left(  -1\right)  ^{n+u}\beta_{u,v}. \label{sol.desnanot.jaw.1}%
\end{align}

Lemma \ref{lem.desnanot.AB.tech} \textbf{(a)} (applied to $w$ instead of $v$)
yields%
\begin{align}
\det\left(  C_{\sim w,\sim q}\right)   &  =-\left(  -1\right)  ^{n+u}%
\underbrace{\det\left(  \operatorname*{rows}\nolimits_{1,2,\ldots
,\widehat{u},\ldots,\widehat{w},\ldots,n}\left(  B_{\bullet,\sim q}\right)
\right)  }_{\substack{=\beta_{u,w}\\\text{(since }\beta_{u,w}=\det\left(
\operatorname*{rows}\nolimits_{1,2,\ldots,\widehat{u},\ldots,\widehat{w}%
,\ldots,n}\left(  B_{\bullet,\sim q}\right)  \right)  \\\text{(by the
definition of }\beta_{u,w}\text{))}}}\nonumber\\
&  =-\left(  -1\right)  ^{n+u}\beta_{u,w}. \label{sol.desnanot.jaw.2}%
\end{align}

Moreover,
\begin{equation}
\operatorname*{sub}\nolimits_{1,2,\ldots,\widehat{v},\ldots,\widehat{w}%
,\ldots,n}^{1,2,\ldots,\widehat{q},\ldots,\widehat{n},\ldots,n}%
C=\operatorname*{rows}\nolimits_{1,2,\ldots,\widehat{v},\ldots,\widehat{w}%
,\ldots,n}\left(  B_{\bullet,\sim q}\right)  \label{sol.desnanot.jaw.3}%
\end{equation}
\footnote{\textit{Proof of (\ref{sol.desnanot.jaw.3}):} From $C=\left(
B\mid\left(  I_{n}\right)  _{\bullet,u}\right)  $, we obtain
\begin{align*}
C_{\bullet,\sim n}  &  =\left(  B\mid\left(  I_{n}\right)  _{\bullet
,u}\right)  _{\bullet,\sim n}=\left(  B\mid\left(  I_{n}\right)  _{\bullet
,u}\right)  _{\bullet,\sim\left(  \left(  n-1\right)  +1\right)
}\ \ \ \ \ \ \ \ \ \ \left(  \text{since }n=\left(  n-1\right)  +1\right) \\
&  =B
\end{align*}
(by Proposition \ref{prop.addcol.props1} \textbf{(d)}, applied to $n-1$, $B$
and $\left(  I_{n}\right)  _{\bullet,u}$ instead of $m$, $A$ and $v$).
\par
Proposition \ref{prop.submatrix.easy} \textbf{(d)} (applied to $n$, $C$,
$n-2$, $\left(  1,2,\ldots,\widehat{v},\ldots,\widehat{w},\ldots,n\right)  $,
$n-2$ and $\left(  1,2,\ldots,\widehat{q},\ldots,\widehat{n},\ldots,n\right)
$ instead of $m$, $A$, $u$, $\left(  i_{1},i_{2},\ldots,i_{u}\right)  $, $v$
and $\left(  j_{1},j_{2},\ldots,j_{v}\right)  $) yields%
\begin{align}
\operatorname*{sub}\nolimits_{1,2,\ldots,\widehat{v},\ldots,\widehat{w}%
,\ldots,n}^{1,2,\ldots,\widehat{q},\ldots,\widehat{n},\ldots,n}C  &
=\operatorname*{rows}\nolimits_{1,2,\ldots,\widehat{v},\ldots,\widehat{w}%
,\ldots,n}\left(  \operatorname*{cols}\nolimits_{1,2,\ldots,\widehat{q}%
,\ldots,\widehat{n},\ldots,n}C\right) \label{sol.desnanot.jaw.3.pf.1}\\
&  =\operatorname*{cols}\nolimits_{1,2,\ldots,\widehat{q},\ldots
,\widehat{n},\ldots,n}\left(  \operatorname*{rows}\nolimits_{1,2,\ldots
,\widehat{v},\ldots,\widehat{w},\ldots,n}C\right)  .\nonumber
\end{align}
\par
But $n$ is a positive integer; hence, $n\in\left\{  1,2,\ldots,n\right\}  $.
Also, $q\in\left\{  1,2,\ldots,n-1\right\}  $. Thus, Proposition
\ref{prop.unrows.basics} \textbf{(h)} (applied to $n$, $C$, $n$ and $q$
instead of $m$, $A$, $v$ and $w$) yields%
\[
\left(  C_{\bullet,\sim n}\right)  _{\bullet,\sim q}=\operatorname*{cols}%
\nolimits_{1,2,\ldots,\widehat{q},\ldots,\widehat{n},\ldots,n}C.
\]
Hence,%
\[
\operatorname*{cols}\nolimits_{1,2,\ldots,\widehat{q},\ldots,\widehat{n}%
,\ldots,n}C=\left(  \underbrace{C_{\bullet,\sim n}}_{=B}\right)
_{\bullet,\sim q}=B_{\bullet,\sim q}.
\]
Hence, (\ref{sol.desnanot.jaw.3.pf.1}) becomes%
\begin{align*}
\operatorname*{sub}\nolimits_{1,2,\ldots,\widehat{v},\ldots,\widehat{w}%
,\ldots,n}^{1,2,\ldots,\widehat{q},\ldots,\widehat{n},\ldots,n}C  &
=\operatorname*{rows}\nolimits_{1,2,\ldots,\widehat{v},\ldots,\widehat{w}%
,\ldots,n}\left(  \underbrace{\operatorname*{cols}\nolimits_{1,2,\ldots
,\widehat{q},\ldots,\widehat{n},\ldots,n}C}_{=B_{\bullet,\sim q}}\right) \\
&  =\operatorname*{rows}\nolimits_{1,2,\ldots,\widehat{v},\ldots
,\widehat{w},\ldots,n}\left(  B_{\bullet,\sim q}\right)  .
\end{align*}
This proves (\ref{sol.desnanot.jaw.3}).}. Taking determinants on both sides of
this equality, we obtain%
\begin{equation}
\det\left(  \operatorname*{sub}\nolimits_{1,2,\ldots,\widehat{v}%
,\ldots,\widehat{w},\ldots,n}^{1,2,\ldots,\widehat{q},\ldots,\widehat{n}%
,\ldots,n}C\right)  =\det\left(  \operatorname*{rows}\nolimits_{1,2,\ldots
,\widehat{v},\ldots,\widehat{w},\ldots,n}\left(  B_{\bullet,\sim q}\right)
\right)  =\beta_{v,w} \label{sol.desnanot.jaw.3b}%
\end{equation}
(since $\beta_{v,w}=\det\left(  \operatorname*{rows}\nolimits_{1,2,\ldots
,\widehat{v},\ldots,\widehat{w},\ldots,n}\left(  B_{\bullet,\sim q}\right)
\right)  $ (by the definition of $\beta_{v,w}$)).

Lemma \ref{lem.desnanot.AB.tech} \textbf{(c)} yields $C_{\sim v,\sim
n}=B_{\sim v,\bullet}$. Thus,%
\begin{equation}
\det\left(  \underbrace{C_{\sim v,\sim n}}_{=B_{\sim v,\bullet}}\right)
=\det\left(  B_{\sim v,\bullet}\right)  =\alpha_{v} \label{sol.desnanot.jaw.4}%
\end{equation}
(since $\alpha_{v}=\det\left(  B_{\sim v,\bullet}\right)  $ (by the definition
of $\alpha_{v}$)).

Lemma \ref{lem.desnanot.AB.tech} \textbf{(c)} (applied to $w$ instead of $v$)
yields $C_{\sim w,\sim n}=B_{\sim w,\bullet}$. Thus,%
\begin{equation}
\det\left(  \underbrace{C_{\sim w,\sim n}}_{=B_{\sim w,\bullet}}\right)
=\det\left(  B_{\sim w,\bullet}\right)  =\alpha_{w} \label{sol.desnanot.jaw.5}%
\end{equation}
(since $\alpha_{w}=\det\left(  B_{\sim w,\bullet}\right)  $ (by the definition
of $\alpha_{w}$)).

Lemma \ref{lem.desnanot.AB.tech} \textbf{(f)} yields%
\begin{equation}
\det C=\left(  -1\right)  ^{n+u}\underbrace{\det\left(  B_{\sim u,\bullet
}\right)  }_{\substack{=\alpha_{u}\\\text{(since }\alpha_{u}=\det\left(
B_{\sim u,\bullet}\right)  \\\text{(by the definition of }\alpha_{u}\text{))}%
}}=\left(  -1\right)  ^{n+u}\alpha_{u}. \label{sol.desnanot.jaw.6}%
\end{equation}

Now, $q\in\left\{  1,2,\ldots,n-1\right\}  $, so that $q\leq n-1<n$. Also,
$n\in\left\{  1,2,\ldots,n\right\}  $ (since $n$ is a positive integer) and
$q\in\left\{  1,2,\ldots,n-1\right\}  \subseteq\left\{  1,2,\ldots,n\right\}
$. Hence, Theorem \ref{thm.desnanot} (applied to $C$, $v$, $w$, $q$ and $n$
instead of $A$, $p$, $q$, $u$ and $v$) yields%
\begin{align*}
&  \det C\cdot\det\left(  \operatorname*{sub}\nolimits_{1,2,\ldots
,\widehat{v},\ldots,\widehat{w},\ldots,n}^{1,2,\ldots,\widehat{q}%
,\ldots,\widehat{n},\ldots,n}C\right) \\
&  =\underbrace{\det\left(  C_{\sim v,\sim q}\right)  }_{\substack{=-\left(
-1\right)  ^{n+u}\beta_{u,v}\\\text{(by (\ref{sol.desnanot.jaw.1}))}}%
}\cdot\underbrace{\det\left(  C_{\sim w,\sim n}\right)  }_{\substack{=\alpha
_{w}\\\text{(by (\ref{sol.desnanot.jaw.5}))}}}-\underbrace{\det\left(  C_{\sim
v,\sim n}\right)  }_{\substack{=\alpha_{v}\\\text{(by
(\ref{sol.desnanot.jaw.4}))}}}\cdot\underbrace{\det\left(  C_{\sim w,\sim
q}\right)  }_{\substack{=-\left(  -1\right)  ^{n+u}\beta_{u,w}\\\text{(by
(\ref{sol.desnanot.jaw.2}))}}}\\
&  =\left(  -\left(  -1\right)  ^{n+u}\beta_{u,v}\right)  \cdot\alpha
_{w}-\alpha_{v}\cdot\left(  -\left(  -1\right)  ^{n+u}\beta_{u,w}\right) \\
&  =-\left(  -1\right)  ^{n+u}\beta_{u,v}\cdot\alpha_{w}+\alpha_{v}%
\cdot\left(  -1\right)  ^{n+u}\beta_{u,w}\\
&  =\alpha_{v}\cdot\left(  -1\right)  ^{n+u}\beta_{u,w}-\left(  -1\right)
^{n+u}\beta_{u,v}\cdot\alpha_{w}=\left(  -1\right)  ^{n+u}\left(  \alpha
_{v}\beta_{u,w}-\beta_{u,v}\alpha_{w}\right)  .
\end{align*}
Hence,%
\begin{align*}
\left(  -1\right)  ^{n+u}\left(  \alpha_{v}\beta_{u,w}-\beta_{u,v}\alpha
_{w}\right)   &  =\underbrace{\det C}_{\substack{=\left(  -1\right)
^{n+u}\alpha_{u}\\\text{(by (\ref{sol.desnanot.jaw.6}))}}}\cdot
\underbrace{\det\left(  \operatorname*{sub}\nolimits_{1,2,\ldots
,\widehat{v},\ldots,\widehat{w},\ldots,n}^{1,2,\ldots,\widehat{q}%
,\ldots,\widehat{n},\ldots,n}C\right)  }_{\substack{=\beta_{v,w}\\\text{(by
(\ref{sol.desnanot.jaw.3b}))}}}\\
&  =\left(  -1\right)  ^{n+u}\alpha_{u}\beta_{v,w}.
\end{align*}
Multiplying both sides of this equality by $\left(  -1\right)  ^{n+u}$, we
find%
\[
\left(  -1\right)  ^{n+u}\left(  -1\right)  ^{n+u}\left(  \alpha_{v}%
\beta_{u,w}-\beta_{u,v}\alpha_{w}\right)  =\underbrace{\left(  -1\right)
^{n+u}\left(  -1\right)  ^{n+u}}_{\substack{=\left(  -1\right)  ^{\left(
n+u\right)  +\left(  n+u\right)  }=1\\\text{(since }\left(  n+u\right)
+\left(  n+u\right)  =2\left(  n+u\right)  \text{ is even)}}}\alpha_{u}%
\beta_{v,w}=\alpha_{u}\beta_{v,w}.
\]
Hence,%
\begin{align*}
\alpha_{u}\beta_{v,w}  &  =\underbrace{\left(  -1\right)  ^{n+u}\left(
-1\right)  ^{n+u}}_{\substack{=\left(  -1\right)  ^{\left(  n+u\right)
+\left(  n+u\right)  }=1\\\text{(since }\left(  n+u\right)  +\left(
n+u\right)  =2\left(  n+u\right)  \text{ is even)}}}\left(  \alpha_{v}%
\beta_{u,w}-\beta_{u,v}\alpha_{w}\right)  =\alpha_{v}\beta_{u,w}-\beta
_{u,v}\alpha_{w}\\
&  =\alpha_{v}\beta_{u,w}-\alpha_{w}\beta_{u,v}.
\end{align*}
Adding $\alpha_{w}\beta_{u,v}$ to both sides of this equality, we obtain
$\alpha_{u}\beta_{v,w}+\alpha_{w}\beta_{u,v}=\alpha_{v}\beta_{u,w}$. This
solves Exercise \ref{exe.desnanot.jaw}.
\end{proof}

\subsection{Solution to Exercise \ref{exe.desnanot.skew}}

Before we solve Exercise \ref{exe.desnanot.skew}], let us state a variant of
Lemma \ref{lem.sol.altern.STAS.2} for the case when $A$ is alternating:

\begin{lemma}
\label{lem.sol.desnanot.skew.STAS=}Let $n\in\mathbb{N}$ and $m\in\mathbb{N}$.
Let $A=\left(  a_{i,j}\right)  _{1\leq i\leq n,\ 1\leq j\leq n}$ be an
alternating $n\times n$-matrix. Let $S=\left(  s_{i,j}\right)  _{1\leq i\leq
n,\ 1\leq j\leq m}$ be an $n\times m$-matrix. Then,%
\[
S^{T}AS=\left(  \sum_{\substack{\left(  k,\ell\right)  \in\left\{
1,2,\ldots,n\right\}  ^{2};\\k<\ell}}\left(  s_{k,i}s_{\ell,j}-s_{\ell
,i}s_{k,j}\right)  a_{k,\ell}\right)  _{1\leq i\leq m,\ 1\leq j\leq m}.
\]

\end{lemma}

\begin{proof}
[Proof of Lemma \ref{lem.sol.desnanot.skew.STAS=}.]Every $\left(  i,j\right)
\in\left\{  1,2,\ldots,m\right\}  ^{2}$ satisfies%
\begin{equation}
\sum_{\left(  k,\ell\right)  \in\left\{  1,2,\ldots,n\right\}  ^{2}}%
s_{k,i}s_{\ell,j}a_{k,\ell}=\sum_{\substack{\left(  k,\ell\right)  \in\left\{
1,2,\ldots,n\right\}  ^{2};\\k<\ell}}\left(  s_{k,i}s_{\ell,j}-s_{\ell
,i}s_{k,j}\right)  a_{k,\ell}. \label{pf.lem.sol.desnanot.skew.STAS=.1}%
\end{equation}

\begin{vershort}
[\textit{Proof of (\ref{pf.lem.sol.desnanot.skew.STAS=.1}):} Let $\left(
i,j\right)  \in\left\{  1,2,\ldots,m\right\}  ^{2}$. Then,%
\begin{align}
\sum_{\substack{\left(  k,\ell\right)  \in\left\{  1,2,\ldots,n\right\}
^{2};\\k=\ell}}s_{k,i}s_{\ell,j}\underbrace{a_{k,\ell}}_{\substack{=a_{\ell
,\ell}\\\text{(since }k=\ell\text{)}}}  &  =\sum_{\substack{\left(
k,\ell\right)  \in\left\{  1,2,\ldots,n\right\}  ^{2};\\k=\ell}}s_{k,i}%
s_{\ell,j}\underbrace{a_{\ell,\ell}}_{\substack{=0\\\text{(by Lemma
\ref{lem.sol.altern.STAS.3} \textbf{(a)}}\\\text{(applied to }i=\ell\text{))}%
}}\nonumber\\
&  =\sum_{\substack{\left(  k,\ell\right)  \in\left\{  1,2,\ldots,n\right\}
^{2};\\k=\ell}}s_{k,i}s_{\ell,j}0=0.
\label{pf.lem.sol.desnanot.skew.STAS=.1.pf.short.1}%
\end{align}
Also,%
\begin{align}
&  \sum_{\substack{\left(  k,\ell\right)  \in\left\{  1,2,\ldots,n\right\}
^{2};\\k>\ell}}s_{k,i}s_{\ell,j}\underbrace{a_{k,\ell}}_{\substack{=-a_{\ell
,k}\\\text{(by Lemma \ref{lem.sol.altern.STAS.3} \textbf{(b)}}\\\text{(applied
to }\left(  i,j\right)  =\left(  k,\ell\right)  \text{))}}}\nonumber\\
&  =-\sum_{\substack{\left(  k,\ell\right)  \in\left\{  1,2,\ldots,n\right\}
^{2};\\k>\ell}}s_{k,i}s_{\ell,j}a_{\ell,k}=-\underbrace{\sum
_{\substack{\left(  \ell,k\right)  \in\left\{  1,2,\ldots,n\right\}
^{2};\\\ell>k}}}_{\substack{=\sum_{\substack{\left(  k,\ell\right)
\in\left\{  1,2,\ldots,n\right\}  ^{2};\\\ell>k}}=\sum_{\substack{\left(
k,\ell\right)  \in\left\{  1,2,\ldots,n\right\}  ^{2};\\k<\ell}}}}s_{\ell
,i}s_{k,j}a_{k,\ell}\nonumber\\
&  \ \ \ \ \ \ \ \ \ \ \left(  \text{here, we have renamed the summation index
}\left(  k,\ell\right)  \text{ as }\left(  \ell,k\right)  \right) \nonumber\\
&  =-\sum_{\substack{\left(  k,\ell\right)  \in\left\{  1,2,\ldots,n\right\}
^{2};\\k<\ell}}s_{\ell,i}s_{k,j}a_{k,\ell}.
\label{pf.lem.sol.desnanot.skew.STAS=.1.pf.short.2}%
\end{align}
Now,%
\begin{align*}
&  \sum_{\left(  k,\ell\right)  \in\left\{  1,2,\ldots,n\right\}  ^{2}}%
s_{k,i}s_{\ell,j}a_{k,\ell}\\
&  =\sum_{\substack{\left(  k,\ell\right)  \in\left\{  1,2,\ldots,n\right\}
^{2};\\k<\ell}}s_{k,i}s_{\ell,j}a_{k,\ell}+\underbrace{\sum_{\substack{\left(
k,\ell\right)  \in\left\{  1,2,\ldots,n\right\}  ^{2};\\k=\ell}}s_{k,i}%
s_{\ell,j}a_{k,\ell}}_{\substack{=0\\\text{(by
(\ref{pf.lem.sol.desnanot.skew.STAS=.1.pf.short.1}))}}}\\
&  \ \ \ \ \ \ \ \ \ \ +\underbrace{\sum_{\substack{\left(  k,\ell\right)
\in\left\{  1,2,\ldots,n\right\}  ^{2};\\k>\ell}}s_{k,i}s_{\ell,j}a_{k,\ell}%
}_{\substack{=-\sum_{\substack{\left(  k,\ell\right)  \in\left\{
1,2,\ldots,n\right\}  ^{2};\\k<\ell}}s_{\ell,i}s_{k,j}a_{k,\ell}\\\text{(by
(\ref{pf.lem.sol.desnanot.skew.STAS=.1.pf.short.2}))}}}\\
&  \ \ \ \ \ \ \ \ \ \ \left(
\begin{array}
[c]{c}%
\text{since each }\left(  k,\ell\right)  \in\left\{  1,2,\ldots,n\right\}
^{2}\text{ satisfies exactly one of}\\
\text{the three assertions }\left(  k<\ell\right)  \text{, }\left(
k=\ell\right)  \text{ and }\left(  k>\ell\right)
\end{array}
\right) \\
&  =\sum_{\substack{\left(  k,\ell\right)  \in\left\{  1,2,\ldots,n\right\}
^{2};\\k<\ell}}s_{k,i}s_{\ell,j}a_{k,\ell}-\sum_{\substack{\left(
k,\ell\right)  \in\left\{  1,2,\ldots,n\right\}  ^{2};\\k<\ell}}s_{\ell
,i}s_{k,j}a_{k,\ell}\\
&  =\sum_{\substack{\left(  k,\ell\right)  \in\left\{  1,2,\ldots,n\right\}
^{2};\\k<\ell}}\left(  s_{k,i}s_{\ell,j}-s_{\ell,i}s_{k,j}\right)  a_{k,\ell}.
\end{align*}
This proves (\ref{pf.lem.sol.desnanot.skew.STAS=.1}).]
\end{vershort}

\begin{verlong}
[\textit{Proof of (\ref{pf.lem.sol.desnanot.skew.STAS=.1}):} Let $\left(
i,j\right)  \in\left\{  1,2,\ldots,m\right\}  ^{2}$. Then,%
\begin{align}
\sum_{\substack{\left(  k,\ell\right)  \in\left\{  1,2,\ldots,n\right\}
^{2};\\k=\ell}}s_{k,i}s_{\ell,j}\underbrace{a_{k,\ell}}_{\substack{=a_{\ell
,\ell}\\\text{(since }k=\ell\text{)}}}  &  =\sum_{\substack{\left(
k,\ell\right)  \in\left\{  1,2,\ldots,n\right\}  ^{2};\\k=\ell}}s_{k,i}%
s_{\ell,j}\underbrace{a_{\ell,\ell}}_{\substack{=0\\\text{(by Lemma
\ref{lem.sol.altern.STAS.3} \textbf{(a)}}\\\text{(applied to }i=\ell\text{))}%
}}\nonumber\\
&  =\sum_{\substack{\left(  k,\ell\right)  \in\left\{  1,2,\ldots,n\right\}
^{2};\\k=\ell}}s_{k,i}s_{\ell,j}0=0.
\label{pf.lem.sol.desnanot.skew.STAS=.1.pf.1}%
\end{align}
Also,%
\begin{align}
&  \sum_{\substack{\left(  k,\ell\right)  \in\left\{  1,2,\ldots,n\right\}
^{2};\\k>\ell}}s_{k,i}s_{\ell,j}\underbrace{a_{k,\ell}}_{\substack{=-a_{\ell
,k}\\\text{(by Lemma \ref{lem.sol.altern.STAS.3} \textbf{(b)}}\\\text{(applied
to }\left(  i,j\right)  =\left(  k,\ell\right)  \text{))}}}\nonumber\\
&  =\sum_{\substack{\left(  k,\ell\right)  \in\left\{  1,2,\ldots,n\right\}
^{2};\\k>\ell}}s_{k,i}s_{\ell,j}\left(  -a_{\ell,k}\right)  =-\sum
_{\substack{\left(  k,\ell\right)  \in\left\{  1,2,\ldots,n\right\}
^{2};\\k>\ell}}s_{k,i}s_{\ell,j}a_{\ell,k}\nonumber\\
&  =-\underbrace{\sum_{\substack{\left(  \ell,k\right)  \in\left\{
1,2,\ldots,n\right\}  ^{2};\\\ell>k}}}_{\substack{=\sum_{\ell\in\left\{
1,2,\ldots,n\right\}  }\sum_{\substack{k\in\left\{  1,2,\ldots,n\right\}
;\\\ell>k}}\\=\sum_{k\in\left\{  1,2,\ldots,n\right\}  }\sum_{\substack{\ell
\in\left\{  1,2,\ldots,n\right\}  ;\\\ell>k}}\\=\sum_{\substack{\left(
k,\ell\right)  \in\left\{  1,2,\ldots,n\right\}  ^{2};\\\ell>k}}=\sum
_{\substack{\left(  k,\ell\right)  \in\left\{  1,2,\ldots,n\right\}
^{2};\\k<\ell}}\\\text{(because for every }\left(  k,\ell\right)  \in\left\{
1,2,\ldots,n\right\}  ^{2}\text{, the}\\\text{statement }\left(
\ell>k\right)  \text{ is equivalent to }\left(  k<\ell\right)  \text{)}%
}}s_{\ell,i}s_{k,j}a_{k,\ell}\nonumber\\
&  \ \ \ \ \ \ \ \ \ \ \left(  \text{here, we have renamed the summation index
}\left(  k,\ell\right)  \text{ as }\left(  \ell,k\right)  \right) \nonumber\\
&  =-\sum_{\substack{\left(  k,\ell\right)  \in\left\{  1,2,\ldots,n\right\}
^{2};\\k<\ell}}s_{\ell,i}s_{k,j}a_{k,\ell}.
\label{pf.lem.sol.desnanot.skew.STAS=.1.pf.2}%
\end{align}
Now,%
\begin{align*}
&  \sum_{\left(  k,\ell\right)  \in\left\{  1,2,\ldots,n\right\}  ^{2}}%
s_{k,i}s_{\ell,j}a_{k,\ell}\\
&  =\sum_{\substack{\left(  k,\ell\right)  \in\left\{  1,2,\ldots,n\right\}
^{2};\\k<\ell}}s_{k,i}s_{\ell,j}a_{k,\ell}+\underbrace{\sum_{\substack{\left(
k,\ell\right)  \in\left\{  1,2,\ldots,n\right\}  ^{2};\\k=\ell}}s_{k,i}%
s_{\ell,j}a_{k,\ell}}_{\substack{=0\\\text{(by
(\ref{pf.lem.sol.desnanot.skew.STAS=.1.pf.1}))}}}\\
&  \ \ \ \ \ \ \ \ \ \ +\underbrace{\sum_{\substack{\left(  k,\ell\right)
\in\left\{  1,2,\ldots,n\right\}  ^{2};\\k>\ell}}s_{k,i}s_{\ell,j}a_{k,\ell}%
}_{\substack{=-\sum_{\substack{\left(  k,\ell\right)  \in\left\{
1,2,\ldots,n\right\}  ^{2};\\k<\ell}}s_{\ell,i}s_{k,j}a_{k,\ell}\\\text{(by
(\ref{pf.lem.sol.desnanot.skew.STAS=.1.pf.2}))}}}\\
&  \ \ \ \ \ \ \ \ \ \ \left(
\begin{array}
[c]{c}%
\text{since each }\left(  k,\ell\right)  \in\left\{  1,2,\ldots,n\right\}
^{2}\text{ satisfies exactly one of}\\
\text{the three assertions }\left(  k<\ell\right)  \text{, }\left(
k=\ell\right)  \text{ and }\left(  k>\ell\right)
\end{array}
\right) \\
&  =\sum_{\substack{\left(  k,\ell\right)  \in\left\{  1,2,\ldots,n\right\}
^{2};\\k<\ell}}s_{k,i}s_{\ell,j}a_{k,\ell}+0+\left(  -\sum_{\substack{\left(
k,\ell\right)  \in\left\{  1,2,\ldots,n\right\}  ^{2};\\k<\ell}}s_{\ell
,i}s_{k,j}a_{k,\ell}\right) \\
&  =\sum_{\substack{\left(  k,\ell\right)  \in\left\{  1,2,\ldots,n\right\}
^{2};\\k<\ell}}s_{k,i}s_{\ell,j}a_{k,\ell}-\sum_{\substack{\left(
k,\ell\right)  \in\left\{  1,2,\ldots,n\right\}  ^{2};\\k<\ell}}s_{\ell
,i}s_{k,j}a_{k,\ell}\\
&  =\sum_{\substack{\left(  k,\ell\right)  \in\left\{  1,2,\ldots,n\right\}
^{2};\\k<\ell}}\left(  s_{k,i}s_{\ell,j}-s_{\ell,i}s_{k,j}\right)  a_{k,\ell}.
\end{align*}
This proves (\ref{pf.lem.sol.desnanot.skew.STAS=.1}).]
\end{verlong}

Lemma \ref{lem.sol.altern.STAS.2} yields%
\begin{align*}
S^{T}AS  &  =\left(  \underbrace{\sum_{\left(  k,\ell\right)  \in\left\{
1,2,\ldots,n\right\}  ^{2}}s_{k,i}s_{\ell,j}a_{k,\ell}}_{\substack{=\sum
_{\substack{\left(  k,\ell\right)  \in\left\{  1,2,\ldots,n\right\}
^{2};\\k<\ell}}\left(  s_{k,i}s_{\ell,j}-s_{\ell,i}s_{k,j}\right)  a_{k,\ell
}\\\text{(by (\ref{pf.lem.sol.desnanot.skew.STAS=.1}))}}}\right)  _{1\leq
i\leq m,\ 1\leq j\leq m}\\
&  =\left(  \sum_{\substack{\left(  k,\ell\right)  \in\left\{  1,2,\ldots
,n\right\}  ^{2};\\k<\ell}}\left(  s_{k,i}s_{\ell,j}-s_{\ell,i}s_{k,j}\right)
a_{k,\ell}\right)  _{1\leq i\leq m,\ 1\leq j\leq m}.
\end{align*}
Thus, Lemma \ref{lem.sol.desnanot.skew.STAS=}.
\end{proof}

Another simple lemma that we will use is the following:

\begin{lemma}
\label{lem.sol.desnanot.skew.multi}Let $n\in\mathbb{N}$. Let $B=\left(
b_{i,j}\right)  _{1\leq i\leq n,\ 1\leq j\leq n}$ be an alternating $n\times
n$-matrix. Let $c\in\mathbb{K}$. Assume that for every $\left(  i,j\right)
\in\left\{  1,2,\ldots,n\right\}  ^{2}$ satisfying $i<j$, the element
$b_{i,j}$ of $\mathbb{K}$ is a multiple of $c$. Then, each entry of $B$ is a
multiple of $c$.
\end{lemma}

\begin{proof}
[Proof of Lemma \ref{lem.sol.desnanot.skew.multi}.]We have assumed that the
following holds:

\begin{statement}
\textit{Fact 1:} For every $\left(  i,j\right)  \in\left\{  1,2,\ldots
,n\right\}  ^{2}$ satisfying $i<j$, the element $b_{i,j}$ of $\mathbb{K}$ is a
multiple of $c$.
\end{statement}

From this, it is easy to derive the following claim:

\begin{statement}
\textit{Fact 2:} For every $\left(  i,j\right)  \in\left\{  1,2,\ldots
,n\right\}  ^{2}$, the element $b_{i,j}$ of $\mathbb{K}$ is a multiple of $c$.
\end{statement}

[\textit{Proof of Fact 2:} Let $\left(  i,j\right)  \in\left\{  1,2,\ldots
,n\right\}  ^{2}$. We must show that the element $b_{i,j}$ of $\mathbb{K}$ is
a multiple of $c$.

We are in one of the following three cases:

\textit{Case 1:} We have $i<j$.

\textit{Case 2:} We have $i=j$.

\textit{Case 3:} We have $i>j$.

Let us first consider Case 1. In this case, we have $i<j$. Hence, the element
$b_{i,j}$ of $\mathbb{K}$ is a multiple of $c$ (by Fact 1). Thus, Fact 2 is
proven in Case 1.

Let us next consider Case 2. In this case, we have $i=j$. Hence, $j=i$.
Therefore, $b_{i,j}=b_{i,i}=0$ (by Lemma \ref{lem.sol.altern.STAS.3}
\textbf{(a)} (applied to $B$ and $b_{u,v}$ instead of $A$ and $a_{u,v}$)).
Thus, the element $b_{i,j}$ of $\mathbb{K}$ is a multiple of $c$ (since the
element $0$ of $\mathbb{K}$ is a multiple of $c$). Hence, Fact 2 is proven in
Case 2.

Let us finally consider Case 3. In this case, we have $i>j$. Thus, $j<i$.
Hence, Fact 1 (applied to $\left(  j,i\right)  $ instead of $\left(
i,j\right)  $) shows that the element $b_{j,i}$ of $\mathbb{K}$ is a multiple
of $c$. In other words, $b_{j,i}=dc$ for some $d\in\mathbb{K}$. Consider this
$d$. But Lemma \ref{lem.sol.altern.STAS.3} \textbf{(b)} (applied to $B$ and
$b_{u,v}$ instead of $A$ and $a_{u,v}$) yields $b_{i,j}=-\underbrace{b_{j,i}%
}_{=dc}=-dc=\left(  -d\right)  c$. Hence, the element $b_{i,j}$ of
$\mathbb{K}$ is a multiple of $c$. Thus, Fact 2 is proven in Case 3.

We have now proven Fact 2 in each of the three Cases 1, 2 and 3. Since these
three Cases cover all possibilities, we thus conclude that Fact 2 always holds.]

Notice that $B=\left(  b_{i,j}\right)  _{1\leq i\leq n,\ 1\leq j\leq n}$.
Thus, for every $\left(  i,j\right)  \in\left\{  1,2,\ldots,n\right\}  ^{2}$,
we have%
\begin{equation}
\left(  \text{the }\left(  i,j\right)  \text{-th entry of }B\right)  =b_{i,j}.
\label{pf.lem.sol.desnanot.skew.multi.1}%
\end{equation}

We now need to prove that each entry of $B$ is a multiple of $c$. In other
words, we need to show that, for every $\left(  i,j\right)  \in\left\{
1,2,\ldots,n\right\}  ^{2}$, the $\left(  i,j\right)  $-th entry of $B$ is a
multiple of $c$. So let us fix $\left(  i,j\right)  \in\left\{  1,2,\ldots
,n\right\}  ^{2}$. Then, Fact 2 shows that $b_{i,j}$ is a multiple of $c$. In
light of (\ref{pf.lem.sol.desnanot.skew.multi.1}), this rewrites as follows:
The $\left(  i,j\right)  $-th entry of $B$ is a multiple of $c$. This is
precisely what we wanted to prove. Thus, the proof of Lemma
\ref{lem.sol.desnanot.skew.multi} is complete.
\end{proof}

\begin{proof}
[Solution to Exercise \ref{exe.desnanot.skew}.]Theorem \ref{thm.desnanot}
leads to the following fact: If $i$, $j$, $k$ and $\ell$ are four elements of
$\left\{  1,2,\ldots,n\right\}  $ such that $i<j$ and $k<\ell$, then%
\begin{align}
&  \det S\cdot\det\left(  \operatorname*{sub}\nolimits_{1,2,\ldots
,\widehat{i},\ldots,\widehat{j},\ldots,n}^{1,2,\ldots,\widehat{k}%
,\ldots,\widehat{\ell},\ldots,n}S\right) \nonumber\\
&  =\det\left(  S_{\sim i,\sim k}\right)  \cdot\det\left(  S_{\sim j,\sim\ell
}\right)  -\det\left(  S_{\sim i,\sim\ell}\right)  \cdot\det\left(  S_{\sim
j,\sim k}\right)  \label{sol.desnanot.skew.des}%
\end{align}
\footnote{\textit{Proof of (\ref{sol.desnanot.skew.des}):} Let $i$, $j$, $k$
and $\ell$ be four elements of $\left\{  1,2,\ldots,n\right\}  $ such that
$i<j$ and $k<\ell$. We have $i\geq1$ (since $i\in\left\{  1,2,\ldots
,n\right\}  $) and $j\leq n$ (since $j\in\left\{  1,2,\ldots,n\right\}  $).
But $i<j$ and thus $i\leq j-1$ (since $i$ and $j$ are integers). Hence, $1\leq
i\leq\underbrace{j}_{\leq n}-1\leq n-1$, so that $n-1\geq1$ and thus $n\geq2$.
Hence, Theorem \ref{thm.desnanot} (applied to $S$, $i$, $j$, $k$ and $\ell$
instead of $A$, $p$, $q$, $u$ and $v$) yields%
\begin{equation}
\det S\cdot\det\left(  \operatorname*{sub}\nolimits_{1,2,\ldots,\widehat{i}%
,\ldots,\widehat{j},\ldots,n}^{1,2,\ldots,\widehat{k},\ldots,\widehat{\ell
},\ldots,n}S\right)  =\det\left(  S_{\sim i,\sim k}\right)  \cdot\det\left(
S_{\sim j,\sim\ell}\right)  -\det\left(  S_{\sim i,\sim\ell}\right)  \cdot
\det\left(  S_{\sim j,\sim k}\right)  .\nonumber
\end{equation}
This proves (\ref{sol.desnanot.skew.des}).}.

For each $\left(  i,j\right)  \in\left\{  1,2,\ldots,n\right\}  ^{2}$, define
an element $s_{i,j}$ of $\mathbb{K}$ by%
\begin{equation}
s_{i,j}=\left(  -1\right)  ^{i+j}\det\left(  S_{\sim j,\sim i}\right)  .
\label{sol.desnanot.skew.1}%
\end{equation}
(Notice that these elements $s_{i,j}$ are \textbf{not} supposed to be the
entries of $S$, despite the notation!)

The definition of $\operatorname*{adj}S$ yields%
\[
\operatorname*{adj}S=\left(  \underbrace{\left(  -1\right)  ^{i+j}\det\left(
S_{\sim j,\sim i}\right)  }_{\substack{=s_{i,j}\\\text{(by
(\ref{sol.desnanot.skew.1}))}}}\right)  _{1\leq i\leq n,\ 1\leq j\leq
n}=\left(  s_{i,j}\right)  _{1\leq i\leq n,\ 1\leq j\leq n}.
\]
Thus, Lemma \ref{lem.sol.desnanot.skew.STAS=} (applied to $n$ and
$\operatorname*{adj}S$ instead of $m$ and $S$) yields%
\begin{equation}
\left(  \operatorname*{adj}S\right)  ^{T}A\left(  \operatorname*{adj}S\right)
=\left(  \sum_{\substack{\left(  k,\ell\right)  \in\left\{  1,2,\ldots
,n\right\}  ^{2};\\k<\ell}}\left(  s_{k,i}s_{\ell,j}-s_{\ell,i}s_{k,j}\right)
a_{k,\ell}\right)  _{1\leq i\leq n,\ 1\leq j\leq n}.
\label{sol.desnanot.skew.3}%
\end{equation}
Moreover, Exercise \ref{exe.altern.STAS} (applied to $n$ and
$\operatorname*{adj}S$ instead of $m$ and $S$) yields that the $n\times
n$-matrix $\left(  \operatorname*{adj}S\right)  ^{T}A\left(
\operatorname*{adj}S\right)  $ is alternating.

But if $i$, $j$, $k$ and $\ell$ are four elements of $\left\{  1,2,\ldots
,n\right\}  $ such that $i<j$ and $k<\ell$, then%
\begin{align}
&  \underbrace{s_{k,i}}_{\substack{=\left(  -1\right)  ^{k+i}\det\left(
S_{\sim i,\sim k}\right)  \\\text{(by the definition of }s_{k,i}\text{)}%
}}\underbrace{s_{\ell,j}}_{\substack{=\left(  -1\right)  ^{\ell+j}\det\left(
S_{\sim j,\sim\ell}\right)  \\\text{(by the definition of }s_{\ell,j}\text{)}%
}}-\underbrace{s_{\ell,i}}_{\substack{=\left(  -1\right)  ^{\ell+i}\det\left(
S_{\sim i,\sim\ell}\right)  \\\text{(by the definition of }s_{\ell,i}\text{)}%
}}\underbrace{s_{k,j}}_{\substack{=\left(  -1\right)  ^{k+j}\det\left(
S_{\sim j,\sim k}\right)  \\\text{(by the definition of }s_{k,j}\text{)}%
}}\nonumber\\
&  =\underbrace{\left(  -1\right)  ^{k+i}\det\left(  S_{\sim i,\sim k}\right)
\cdot\left(  -1\right)  ^{\ell+j}\det\left(  S_{\sim j,\sim\ell}\right)
}_{=\left(  -1\right)  ^{k+i}\left(  -1\right)  ^{\ell+j}\det\left(  S_{\sim
i,\sim k}\right)  \cdot\det\left(  S_{\sim j,\sim\ell}\right)  }%
-\underbrace{\left(  -1\right)  ^{\ell+i}\det\left(  S_{\sim i,\sim\ell
}\right)  \cdot\left(  -1\right)  ^{k+j}\det\left(  S_{\sim j,\sim k}\right)
}_{=\left(  -1\right)  ^{\ell+i}\left(  -1\right)  ^{k+j}\det\left(  S_{\sim
i,\sim\ell}\right)  \cdot\det\left(  S_{\sim j,\sim k}\right)  }\nonumber\\
&  =\underbrace{\left(  -1\right)  ^{k+i}\left(  -1\right)  ^{\ell+j}%
}_{\substack{=\left(  -1\right)  ^{\left(  k+i\right)  +\left(  \ell+j\right)
}=\left(  -1\right)  ^{k+\ell+i+j}\\\text{(since }\left(  k+i\right)  +\left(
\ell+j\right)  =k+\ell+i+j\text{)}}}\det\left(  S_{\sim i,\sim k}\right)
\cdot\det\left(  S_{\sim j,\sim\ell}\right) \nonumber\\
&  \ \ \ \ \ \ \ \ \ \ -\underbrace{\left(  -1\right)  ^{\ell+i}\left(
-1\right)  ^{k+j}}_{\substack{=\left(  -1\right)  ^{\left(  \ell+i\right)
+\left(  k+j\right)  }=\left(  -1\right)  ^{k+\ell+i+j}\\\text{(since }\left(
\ell+i\right)  +\left(  k+j\right)  =k+\ell+i+j\text{)}}}\det\left(  S_{\sim
i,\sim\ell}\right)  \cdot\det\left(  S_{\sim j,\sim k}\right) \nonumber\\
&  =\left(  -1\right)  ^{k+\ell+i+j}\det\left(  S_{\sim i,\sim k}\right)
\cdot\det\left(  S_{\sim j,\sim\ell}\right)  -\left(  -1\right)  ^{k+\ell
+i+j}\det\left(  S_{\sim i,\sim\ell}\right)  \cdot\det\left(  S_{\sim j,\sim
k}\right) \nonumber\\
&  =\left(  -1\right)  ^{k+\ell+i+j}\underbrace{\left(  \det\left(  S_{\sim
i,\sim k}\right)  \cdot\det\left(  S_{\sim j,\sim\ell}\right)  -\det\left(
S_{\sim i,\sim\ell}\right)  \cdot\det\left(  S_{\sim j,\sim k}\right)
\right)  }_{\substack{=\det S\cdot\det\left(  \operatorname*{sub}%
\nolimits_{1,2,\ldots,\widehat{i},\ldots,\widehat{j},\ldots,n}^{1,2,\ldots
,\widehat{k},\ldots,\widehat{\ell},\ldots,n}S\right)  \\\text{(by
(\ref{sol.desnanot.skew.des}))}}}\nonumber\\
&  =\left(  -1\right)  ^{k+\ell+i+j}\det S\cdot\det\left(  \operatorname*{sub}%
\nolimits_{1,2,\ldots,\widehat{i},\ldots,\widehat{j},\ldots,n}^{1,2,\ldots
,\widehat{k},\ldots,\widehat{\ell},\ldots,n}S\right)  .
\label{sol.desnanot.skew.5}%
\end{align}
Hence, for every $\left(  i,j\right)  \in\left\{  1,2,\ldots,n\right\}  ^{2}$
satisfying $i<j$, the element \newline$\sum_{\substack{\left(  k,\ell\right)
\in\left\{  1,2,\ldots,n\right\}  ^{2};\\k<\ell}}\left(  s_{k,i}s_{\ell
,j}-s_{\ell,i}s_{k,j}\right)  a_{k,\ell}$ of $\mathbb{K}$ is a multiple of
$\det S$\ \ \ \ \footnote{\textit{Proof.} Let $\left(  i,j\right)  \in\left\{
1,2,\ldots,n\right\}  ^{2}$ be such that $i<j$. From $\left(  i,j\right)
\in\left\{  1,2,\ldots,n\right\}  ^{2}$, we obtain $i\in\left\{
1,2,\ldots,n\right\}  $ and $j\in\left\{  1,2,\ldots,n\right\}  $. Now,%
\begin{align*}
&  \sum_{\substack{\left(  k,\ell\right)  \in\left\{  1,2,\ldots,n\right\}
^{2};\\k<\ell}}\underbrace{\left(  s_{k,i}s_{\ell,j}-s_{\ell,i}s_{k,j}\right)
}_{\substack{=\left(  -1\right)  ^{k+\ell+i+j}\det S\cdot\det\left(
\operatorname*{sub}\nolimits_{1,2,\ldots,\widehat{i},\ldots,\widehat{j}%
,\ldots,n}^{1,2,\ldots,\widehat{k},\ldots,\widehat{\ell},\ldots,n}S\right)
\\\text{(by (\ref{sol.desnanot.skew.5}))}}}a_{k,\ell}\\
&  =\sum_{\substack{\left(  k,\ell\right)  \in\left\{  1,2,\ldots,n\right\}
^{2};\\k<\ell}}\left(  -1\right)  ^{k+\ell+i+j}\det S\cdot\det\left(
\operatorname*{sub}\nolimits_{1,2,\ldots,\widehat{i},\ldots,\widehat{j}%
,\ldots,n}^{1,2,\ldots,\widehat{k},\ldots,\widehat{\ell},\ldots,n}S\right)
a_{k,\ell}\\
&  =\left(  \sum_{\substack{\left(  k,\ell\right)  \in\left\{  1,2,\ldots
,n\right\}  ^{2};\\k<\ell}}\left(  -1\right)  ^{k+\ell+i+j}\det\left(
\operatorname*{sub}\nolimits_{1,2,\ldots,\widehat{i},\ldots,\widehat{j}%
,\ldots,n}^{1,2,\ldots,\widehat{k},\ldots,\widehat{\ell},\ldots,n}S\right)
a_{k,\ell}\right)  \det S
\end{align*}
is clearly a multiple of $\det S$. Qed.}. Thus, Lemma
\ref{lem.sol.desnanot.skew.multi} (applied to $B=\left(  \operatorname*{adj}%
S\right)  ^{T}A\left(  \operatorname*{adj}S\right)  $, $b_{i,j}=\sum
_{\substack{\left(  k,\ell\right)  \in\left\{  1,2,\ldots,n\right\}
^{2};\\k<\ell}}\left(  s_{k,i}s_{\ell,j}-s_{\ell,i}s_{k,j}\right)  a_{k,\ell}$
and $c=\det S$) shows that each entry of the matrix $\left(
\operatorname*{adj}S\right)  ^{T}A\left(  \operatorname*{adj}S\right)  $ is a
multiple of $\det S$ (because the matrix $\left(  \operatorname*{adj}S\right)
^{T}A\left(  \operatorname*{adj}S\right)  $ is alternating and satisfies
(\ref{sol.desnanot.skew.3})). This solves Exercise \ref{exe.desnanot.skew}.
\end{proof}

\subsection{Solution to Exercise \ref{exe.pluecker.rederive-AC}}

\begin{proof}
[Second proof of Proposition \ref{prop.desnanot.AC}.]Recall that $I_{m}$
denotes the $m\times m$ identity matrix for each $m\in\mathbb{N}$. Thus,
$I_{n}$ is the $n\times n$ identity matrix. Hence, $\left(  I_{n}\right)
_{\bullet,v}$ is a well-defined $n\times1$-matrix (since $v\in\left\{
1,2,\ldots,n\right\}  $).

We know that $C$ is an $n\times n$-matrix (since $C\in\mathbb{K}^{n\times n}%
$), and that $\left(  I_{n}\right)  _{\bullet,v}$ is an $n\times1$-matrix.
Thus, $\left(  C\mid\left(  I_{n}\right)  _{\bullet,v}\right)  $ is an
$n\times\left(  n+1\right)  $-matrix. In other words, $\left(  C\mid\left(
I_{n}\right)  _{\bullet,v}\right)  \in\mathbb{K}^{n\times\left(  n+1\right)
}$.

Set $B=\left(  C\mid\left(  I_{n}\right)  _{\bullet,v}\right)  $. Thus,
$B=\left(  C\mid\left(  I_{n}\right)  _{\bullet,v}\right)  \in\mathbb{K}%
^{n\times\left(  n+1\right)  }$. We have%
\begin{equation}
B_{\bullet,r}=C_{\bullet,r}\ \ \ \ \ \ \ \ \ \ \text{for every }r\in\left\{
1,2,\ldots,n\right\}  \label{sol.pluecker.rederive-AC.1}%
\end{equation}
\footnote{\textit{Proof of (\ref{sol.pluecker.rederive-AC.1}):} Let
$r\in\left\{  1,2,\ldots,n\right\}  $. Proposition \ref{prop.addcol.props1}
\textbf{(a)} (applied to $n$, $C$, $\left(  I_{n}\right)  _{\bullet,v}$ and
$r$ instead of $m$, $A$, $v$ and $q$) yields $\left(  C\mid\left(
I_{n}\right)  _{\bullet,v}\right)  _{\bullet,r}=C_{\bullet,r}$. Now, from
$B=\left(  C\mid\left(  I_{n}\right)  _{\bullet,v}\right)  $, we obtain
$B_{\bullet,r}=\left(  C\mid\left(  I_{n}\right)  _{\bullet,v}\right)
_{\bullet,r}=C_{\bullet,r}$. This proves (\ref{sol.pluecker.rederive-AC.1}).}.
Furthermore,%
\begin{equation}
\det\left(  A\mid B_{\bullet,n+1}\right)  =\left(  -1\right)  ^{n+v}%
\det\left(  A_{\sim v,\bullet}\right)  \label{sol.pluecker.rederive-AC.2}%
\end{equation}
\footnote{\textit{Proof of (\ref{sol.pluecker.rederive-AC.2}):} Proposition
\ref{prop.addcol.props1} \textbf{(b)} (applied to $n$, $C$ and $\left(
I_{n}\right)  _{\bullet,v}$ instead of $m$, $A$ and $v$) yields $\left(
C\mid\left(  I_{n}\right)  _{\bullet,v}\right)  _{\bullet,n+1}=\left(
I_{n}\right)  _{\bullet,v}$. Now, from $B=\left(  C\mid\left(  I_{n}\right)
_{\bullet,v}\right)  $, we obtain $B_{\bullet,n+1}=\left(  C\mid\left(
I_{n}\right)  _{\bullet,v}\right)  _{\bullet,n+1}=\left(  I_{n}\right)
_{\bullet,v}$. Hence,%
\[
\det\left(  A\mid\underbrace{B_{\bullet,n+1}}_{=\left(  I_{n}\right)
_{\bullet,v}}\right)  =\det\left(  A\mid\left(  I_{n}\right)  _{\bullet
,v}\right)  =\left(  -1\right)  ^{n+v}\det\left(  A_{\sim v,\bullet}\right)
\]
(by Proposition \ref{prop.addcol.props2} \textbf{(b)}, applied to $p=v$). This
proves (\ref{sol.pluecker.rederive-AC.2}).}. Moreover,%
\begin{equation}
\det\left(  B_{\bullet,\sim r}\right)  =\left(  -1\right)  ^{n+v}\det\left(
C_{\sim v,\sim r}\right)  \ \ \ \ \ \ \ \ \ \ \text{for every }r\in\left\{
1,2,\ldots,n\right\}  \label{sol.pluecker.rederive-AC.3}%
\end{equation}
\footnote{\textit{Proof of (\ref{sol.pluecker.rederive-AC.3}):} Let
$r\in\left\{  1,2,\ldots,n\right\}  $. From $B=\left(  C\mid\left(
I_{n}\right)  _{\bullet,v}\right)  $, we obtain $B_{\bullet,\sim r}=\left(
C\mid\left(  I_{n}\right)  _{\bullet,v}\right)  _{\bullet,\sim r}=\left(
C_{\bullet,\sim r}\mid\left(  I_{n}\right)  _{\bullet,v}\right)  $ (by
Proposition \ref{prop.addcol.props1} \textbf{(c)}, applied to $n$, $C$,
$\left(  I_{n}\right)  _{\bullet,v}$ and $r$ instead of $m$, $A$, $v$ and
$q$). Hence,%
\begin{equation}
\det\left(  \underbrace{B_{\bullet,\sim r}}_{=\left(  C_{\bullet,\sim r}%
\mid\left(  I_{n}\right)  _{\bullet,v}\right)  }\right)  =\det\left(
C_{\bullet,\sim r}\mid\left(  I_{n}\right)  _{\bullet,v}\right)  =\left(
-1\right)  ^{n+v}\det\left(  \left(  C_{\bullet,\sim r}\right)  _{\sim
v,\bullet}\right)  \label{sol.pluecker.rederive-AC.3.pf.1}%
\end{equation}
(by Proposition \ref{prop.addcol.props2} \textbf{(b)}, applied to $v$ and
$C_{\bullet,\sim r}$ instead of $p$ and $A$). But Proposition
\ref{prop.unrows.basics} \textbf{(c)} (applied to $n$, $C$, $v$ and $r$
instead of $m$, $A$, $u$ and $v$) yields $\left(  C_{\bullet,\sim r}\right)
_{\sim v,\bullet}=\left(  C_{\sim v,\bullet}\right)  _{\bullet,\sim r}=C_{\sim
v,\sim r}$. Now, (\ref{sol.pluecker.rederive-AC.3.pf.1}) becomes%
\[
\det\left(  B_{\bullet,\sim r}\right)  =\left(  -1\right)  ^{n+v}\det\left(
\underbrace{\left(  C_{\bullet,\sim r}\right)  _{\sim v,\bullet}}_{=C_{\sim
v,\sim r}}\right)  =\left(  -1\right)  ^{n+v}\det\left(  C_{\sim v,\sim
r}\right)  .
\]
This proves (\ref{sol.pluecker.rederive-AC.3}).}. Finally,%
\begin{equation}
B_{\bullet,\sim\left(  n+1\right)  }=C \label{sol.pluecker.rederive-AC.4}%
\end{equation}
\footnote{\textit{Proof of (\ref{sol.pluecker.rederive-AC.4}):} Proposition
\ref{prop.addcol.props1} \textbf{(d)} (applied to $n$, $C$ and $\left(
I_{n}\right)  _{\bullet,v}$ instead of $m$, $A$ and $v$) yields $\left(
C\mid\left(  I_{n}\right)  _{\bullet,v}\right)  _{\bullet,\sim\left(
n+1\right)  }=C$. Now, from $B=\left(  C\mid\left(  I_{n}\right)  _{\bullet
,v}\right)  $, we obtain $B_{\bullet,\sim\left(  n+1\right)  }=\left(
C\mid\left(  I_{n}\right)  _{\bullet,v}\right)  _{\bullet,\sim\left(
n+1\right)  }=C$. This proves (\ref{sol.pluecker.rederive-AC.4}).}.

Now, Theorem \ref{thm.pluecker.plu} yields%
\begin{align*}
0  &  =\sum_{r=1}^{n+1}\left(  -1\right)  ^{r}\det\left(  A\mid B_{\bullet
,r}\right)  \det\left(  B_{\bullet,\sim r}\right) \\
&  =\sum_{r=1}^{n}\left(  -1\right)  ^{r}\det\left(  A\mid
\underbrace{B_{\bullet,r}}_{\substack{=C_{\bullet,r}\\\text{(by
(\ref{sol.pluecker.rederive-AC.1}))}}}\right)  \underbrace{\det\left(
B_{\bullet,\sim r}\right)  }_{\substack{=\left(  -1\right)  ^{n+v}\det\left(
C_{\sim v,\sim r}\right)  \\\text{(by (\ref{sol.pluecker.rederive-AC.3}))}}}\\
&  \ \ \ \ \ \ \ \ \ \ +\left(  -1\right)  ^{n+1}\underbrace{\det\left(  A\mid
B_{\bullet,n+1}\right)  }_{\substack{=\left(  -1\right)  ^{n+v}\det\left(
A_{\sim v,\bullet}\right)  \\\text{(by (\ref{sol.pluecker.rederive-AC.2}))}%
}}\det\left(  \underbrace{B_{\bullet,\sim\left(  n+1\right)  }}%
_{\substack{=C\\\text{(by (\ref{sol.pluecker.rederive-AC.4}))}}}\right) \\
&  \ \ \ \ \ \ \ \ \ \ \left(  \text{here, we have split off the addend for
}r=n+1\text{ from the sum}\right) \\
&  =\sum_{r=1}^{n}\left(  -1\right)  ^{r}\underbrace{\det\left(  A\mid
C_{\bullet,r}\right)  \left(  -1\right)  ^{n+v}}_{=\left(  -1\right)
^{n+v}\det\left(  A\mid C_{\bullet,r}\right)  }\det\left(  C_{\sim v,\sim
r}\right) \\
&  \ \ \ \ \ \ \ \ \ \ +\underbrace{\left(  -1\right)  ^{n+1}\left(
-1\right)  ^{n+v}}_{\substack{=\left(  -1\right)  ^{\left(  n+1\right)
+\left(  n+v\right)  }=\left(  -1\right)  ^{v+1}\\\text{(since }\left(
n+1\right)  +\left(  n+v\right)  =2n+\left(  v+1\right)  \equiv
v+1\operatorname{mod}2\text{)}}}\det\left(  A_{\sim v,\bullet}\right)  \det
C\\
&  =\sum_{r=1}^{n}\underbrace{\left(  -1\right)  ^{r}\left(  -1\right)
^{n+v}}_{=\left(  -1\right)  ^{r+\left(  n+v\right)  }}\det\left(  A\mid
C_{\bullet,r}\right)  \det\left(  C_{\sim v,\sim r}\right)
+\underbrace{\left(  -1\right)  ^{v+1}}_{=-\left(  -1\right)  ^{v}}\det\left(
A_{\sim v,\bullet}\right)  \det C\\
&  =\sum_{r=1}^{n}\left(  -1\right)  ^{r+\left(  n+v\right)  }\det\left(
A\mid C_{\bullet,r}\right)  \det\left(  C_{\sim v,\sim r}\right)  -\left(
-1\right)  ^{v}\det\left(  A_{\sim v,\bullet}\right)  \det C.
\end{align*}
Adding $\left(  -1\right)  ^{v}\det\left(  A_{\sim v,\bullet}\right)  \det C$
to both sides of this equality, we obtain%
\[
\left(  -1\right)  ^{v}\det\left(  A_{\sim v,\bullet}\right)  \det
C=\sum_{r=1}^{n}\left(  -1\right)  ^{r+\left(  n+v\right)  }\det\left(  A\mid
C_{\bullet,r}\right)  \det\left(  C_{\sim v,\sim r}\right)  .
\]
Multiplying both sides of this equality by $\left(  -1\right)  ^{v}$, we find%
\begin{align*}
&  \left(  -1\right)  ^{v}\left(  -1\right)  ^{v}\det\left(  A_{\sim
v,\bullet}\right)  \det C\\
&  =\left(  -1\right)  ^{v}\sum_{r=1}^{n}\left(  -1\right)  ^{r+\left(
n+v\right)  }\det\left(  A\mid C_{\bullet,r}\right)  \det\left(  C_{\sim
v,\sim r}\right) \\
&  =\sum_{r=1}^{n}\underbrace{\left(  -1\right)  ^{v}\left(  -1\right)
^{r+\left(  n+v\right)  }}_{\substack{=\left(  -1\right)  ^{v+\left(
r+\left(  n+v\right)  \right)  }=\left(  -1\right)  ^{n+r}\\\text{(since
}v+\left(  r+\left(  n+v\right)  \right)  =2v+n+r\equiv n+r\operatorname{mod}%
2\text{)}}}\det\left(  A\mid C_{\bullet,r}\right)  \det\left(  C_{\sim v,\sim
r}\right) \\
&  =\sum_{r=1}^{n}\left(  -1\right)  ^{n+r}\det\left(  A\mid C_{\bullet
,r}\right)  \det\left(  C_{\sim v,\sim r}\right) \\
&  =\sum_{q=1}^{n}\left(  -1\right)  ^{n+q}\det\left(  A\mid C_{\bullet
,q}\right)  \det\left(  C_{\sim v,\sim q}\right)
\end{align*}
(here, we have renamed the summation index $r$ as $q$). Hence,%
\begin{align*}
&  \sum_{q=1}^{n}\left(  -1\right)  ^{n+q}\det\left(  A\mid C_{\bullet
,q}\right)  \det\left(  C_{\sim v,\sim q}\right) \\
&  =\underbrace{\left(  -1\right)  ^{v}\left(  -1\right)  ^{v}}%
_{\substack{=\left(  -1\right)  ^{v+v}=1\\\text{(since }v+v=2v\text{ is
even)}}}\det\left(  A_{\sim v,\bullet}\right)  \det C=\det\left(  A_{\sim
v,\bullet}\right)  \det C.
\end{align*}
This proves Proposition \ref{prop.desnanot.AC} again.
\end{proof}

Thus, Exercise \ref{exe.pluecker.rederive-AC} is solved.

\subsection{Solution to Exercise \ref{exe.det.laplace-multi}}

Throughout this section, we shall use the notations introduced in Definition
\ref{def.submatrix} and in Definition \ref{def.sect.laplace.notations}.

Let us prepare for the proof of Lemma \ref{lem.det.laplace-multi.Apq} by
showing a particular case:

\begin{lemma}
\label{lem.sol.det.laplace-multi.1}Let $n\in\mathbb{N}$. For any subset $I$ of
$\left\{  1,2,\ldots,n\right\}  $, we let $\widetilde{I}$ denote the
complement $\left\{  1,2,\ldots,n\right\}  \setminus I$ of $I$.

Let $A=\left(  a_{i,j}\right)  _{1\leq i\leq n,\ 1\leq j\leq n}$ and
$B=\left(  b_{i,j}\right)  _{1\leq i\leq n,\ 1\leq j\leq n}$ be two $n\times
n$-matrices. Let $Q$ be a subset of $\left\{  1,2,\ldots,n\right\}  $. Let
$k=\left\vert Q\right\vert $. Then,%
\begin{align*}
&  \sum_{\substack{\sigma\in S_{n};\\\sigma\left(  \left\{  1,2,\ldots
,k\right\}  \right)  =Q}}\left(  -1\right)  ^{\sigma}\left(  \prod
_{i\in\left\{  1,2,\ldots,k\right\}  }a_{i,\sigma\left(  i\right)  }\right)
\left(  \prod_{i\in\left\{  k+1,k+2,\ldots,n\right\}  }b_{i,\sigma\left(
i\right)  }\right) \\
&  =\left(  -1\right)  ^{\left(  1+2+\cdots+k\right)  +\sum Q}\det\left(
\operatorname*{sub}\nolimits_{\left(  1,2,\ldots,k\right)  }^{w\left(
Q\right)  }A\right)  \det\left(  \operatorname*{sub}\nolimits_{\left(
k+1,k+2,\ldots,n\right)  }^{w\left(  \widetilde{Q}\right)  }B\right)  .
\end{align*}

\end{lemma}

\begin{proof}
[Proof of Lemma \ref{lem.sol.det.laplace-multi.1}.]We begin with some simple observations.

\begin{vershort}
The definition of $\widetilde{Q}$ yields $\widetilde{Q}=\left\{
1,2,\ldots,n\right\}  \setminus Q$. Since $Q$ is a subset of $\left\{
1,2,\ldots,n\right\}  $, this leads to $\left\vert \widetilde{Q}\right\vert
=\underbrace{\left\vert \left\{  1,2,\ldots,n\right\}  \right\vert }%
_{=n}-\underbrace{\left\vert Q\right\vert }_{=k}=n-k$. Thus, $n-k=\left\vert
\widetilde{Q}\right\vert \geq0$, so that $n\geq k$ and thus $k\in\left\{
0,1,\ldots,n\right\}  $.
\end{vershort}

\begin{verlong}
The definition of $\widetilde{Q}$ yields $\widetilde{Q}=\left\{
1,2,\ldots,n\right\}  \setminus Q\subseteq\left\{  1,2,\ldots,n\right\}  $.
Thus, $\widetilde{Q}$ is a finite set (since $\left\{  1,2,\ldots,n\right\}  $
is a finite set). Hence, $\left\vert \widetilde{Q}\right\vert \in\mathbb{N}$.

Also, $Q\subseteq\left\{  1,2,\ldots,n\right\}  $. Thus, $Q$ is a finite set
(since $\left\{  1,2,\ldots,n\right\}  $ is a finite set). Hence, $\left\vert
Q\right\vert \in\mathbb{N}$.

Also,%
\begin{align*}
\left\vert \underbrace{\widetilde{Q}}_{=\left\{  1,2,\ldots,n\right\}
\setminus Q}\right\vert  &  =\left\vert \left\{  1,2,\ldots,n\right\}
\setminus Q\right\vert =\underbrace{\left\vert \left\{  1,2,\ldots,n\right\}
\right\vert }_{=n}-\underbrace{\left\vert Q\right\vert }_{=k}%
\ \ \ \ \ \ \ \ \ \ \left(  \text{since }Q\subseteq\left\{  1,2,\ldots
,n\right\}  \right) \\
&  =n-k.
\end{align*}
Hence, $n-k=\left\vert \widetilde{Q}\right\vert \in\mathbb{N}$; thus,
$n-k\geq0$. In other words, $k\leq n$. Combined with $k\geq0$ (since
$k=\left\vert Q\right\vert \in\mathbb{N}$), this yields $k\in\left\{
0,1,\ldots,n\right\}  $.
\end{verlong}

Let $\left(  q_{1},q_{2},\ldots,q_{k}\right)  $ be the list of all elements of
$Q$ in increasing order (with no repetitions). (This is well-defined (by
Definition \ref{def.ind.inclist}), because $\left\vert Q\right\vert =k$.)

Let $\left(  r_{1},r_{2},\ldots,r_{n-k}\right)  $ be the list of all elements
of $\widetilde{Q}$ in increasing order (with no repetitions). (This is
well-defined, because $\left\vert \widetilde{Q}\right\vert =n-k$.)

\begin{vershort}
The lists $w\left(  Q\right)  $ and $\left(  q_{1},q_{2},\ldots,q_{k}\right)
$ must be identical (since they are both defined to be the list of all
elements of $Q$ in increasing order (with no repetitions)). In other words, we
have $w\left(  Q\right)  =\left(  q_{1},q_{2},\ldots,q_{k}\right)  $.
Similarly, $w\left(  \widetilde{Q}\right)  =\left(  r_{1},r_{2},\ldots
,r_{n-k}\right)  $.
\end{vershort}

\begin{verlong}
But $w\left(  Q\right)  $ is the list of all elements of $Q$ in increasing
order (with no repetitions) (by the definition of $w\left(  Q\right)  $).
Thus,%
\begin{align*}
w\left(  Q\right)   &  =\left(  \text{the list of all elements of }Q\text{ in
increasing order (with no repetitions)}\right) \\
&  =\left(  q_{1},q_{2},\ldots,q_{k}\right)
\end{align*}
(since $\left(  q_{1},q_{2},\ldots,q_{k}\right)  $ is the list of all elements
of $Q$ in increasing order (with no repetitions)).

Also, $w\left(  \widetilde{Q}\right)  $ is the list of all elements of
$\widetilde{Q}$ in increasing order (with no repetitions) (by the definition
of $w\left(  \widetilde{Q}\right)  $). Thus,%
\begin{align*}
w\left(  \widetilde{Q}\right)   &  =\left(  \text{the list of all elements of
}\widetilde{Q}\text{ in increasing order (with no repetitions)}\right) \\
&  =\left(  r_{1},r_{2},\ldots,r_{n-k}\right)
\end{align*}
(since $\left(  r_{1},r_{2},\ldots,r_{n-k}\right)  $ is the list of all
elements of $\widetilde{Q}$ in increasing order (with no repetitions)).
\end{verlong}

We know that $\left(  r_{1},r_{2},\ldots,r_{n-k}\right)  $ is the list of all
elements of $\widetilde{Q}$ in increasing order (with no repetitions). In
other words, $\left(  r_{1},r_{2},\ldots,r_{n-k}\right)  $ is the list of all
elements of $\left\{  1,2,\ldots,n\right\}  \setminus Q$ in increasing order
(with no repetitions) (since $\widetilde{Q}=\left\{  1,2,\ldots,n\right\}
\setminus Q$). Thus, for every $\alpha\in S_{k}$ and $\beta\in S_{n-k}$, we
can define an element $\sigma_{Q,\alpha,\beta}\in S_{n}$ according to Exercise
\ref{exe.Ialbe} \textbf{(a)} (applied to $Q$, $\left(  q_{1},q_{2}%
,\ldots,q_{k}\right)  $ and $\left(  r_{1},r_{2},\ldots,r_{n-k}\right)  $
instead of $I$, $\left(  a_{1},a_{2},\ldots,a_{k}\right)  $ and $\left(
b_{1},b_{2},\ldots,b_{n-k}\right)  $). Consider this $\sigma_{Q,\alpha,\beta}%
$. Exercise \ref{exe.Ialbe} \textbf{(b)} (applied to $Q$, $\left(  q_{1}%
,q_{2},\ldots,q_{k}\right)  $ and $\left(  r_{1},r_{2},\ldots,r_{n-k}\right)
$ instead of $I$, $\left(  a_{1},a_{2},\ldots,a_{k}\right)  $ and $\left(
b_{1},b_{2},\ldots,b_{n-k}\right)  $) shows that for every $\alpha\in S_{k}$
and $\beta\in S_{n-k}$, we have%
\[
\ell\left(  \sigma_{Q,\alpha,\beta}\right)  =\ell\left(  \alpha\right)
+\ell\left(  \beta\right)  +\sum Q-\left(  1+2+\cdots+k\right)
\]
and%
\begin{equation}
\left(  -1\right)  ^{\sigma_{Q,\alpha,\beta}}=\left(  -1\right)  ^{\alpha
}\cdot\left(  -1\right)  ^{\beta}\cdot\left(  -1\right)  ^{\sum Q-\left(
1+2+\cdots+k\right)  }. \label{pf.lem.sol.det.laplace-multi.1.sign}%
\end{equation}

Exercise \ref{exe.Ialbe} \textbf{(c)} (applied to $Q$, $\left(  q_{1}%
,q_{2},\ldots,q_{k}\right)  $ and $\left(  r_{1},r_{2},\ldots,r_{n-k}\right)
$ instead of $I$, $\left(  a_{1},a_{2},\ldots,a_{k}\right)  $ and $\left(
b_{1},b_{2},\ldots,b_{n-k}\right)  $) shows that the map%
\begin{align*}
S_{k}\times S_{n-k}  &  \rightarrow\left\{  \tau\in S_{n}\ \mid\ \tau\left(
\left\{  1,2,\ldots,k\right\}  \right)  =Q\right\}  ,\\
\left(  \alpha,\beta\right)   &  \mapsto\sigma_{Q,\alpha,\beta}%
\end{align*}
is well-defined and a bijection.

We have $w\left(  Q\right)  =\left(  q_{1},q_{2},\ldots,q_{k}\right)  $. Thus,%
\begin{align*}
\operatorname*{sub}\nolimits_{\left(  1,2,\ldots,k\right)  }^{w\left(
Q\right)  }A  &  =\operatorname*{sub}\nolimits_{\left(  1,2,\ldots,k\right)
}^{\left(  q_{1},q_{2},\ldots,q_{k}\right)  }A=\operatorname*{sub}%
\nolimits_{1,2,\ldots,k}^{q_{1},q_{2},\ldots,q_{k}}A=\left(  a_{x,q_{y}%
}\right)  _{1\leq x\leq k,\ 1\leq y\leq k}\\
&  \ \ \ \ \ \ \ \ \ \ \left(
\begin{array}
[c]{c}%
\text{by the definition of }\operatorname*{sub}\nolimits_{1,2,\ldots,k}%
^{q_{1},q_{2},\ldots,q_{k}}A\text{,}\\
\text{since }A=\left(  a_{i,j}\right)  _{1\leq i\leq n,\ 1\leq j\leq n}%
\end{array}
\right) \\
&  =\left(  a_{i,q_{j}}\right)  _{1\leq i\leq k,\ 1\leq j\leq k}%
\end{align*}
(here, we have renamed the index $\left(  x,y\right)  $ as $\left(
i,j\right)  $). Thus,%
\begin{align}
\det\left(  \operatorname*{sub}\nolimits_{\left(  1,2,\ldots,k\right)
}^{w\left(  Q\right)  }A\right)   &  =\sum_{\sigma\in S_{k}}\left(  -1\right)
^{\sigma}\prod_{i=1}^{k}a_{i,q_{\sigma\left(  i\right)  }}\nonumber\\
&  \ \ \ \ \ \ \ \ \ \ \left(
\begin{array}
[c]{c}%
\text{by (\ref{eq.det.eq.2}), applied to }k\text{, }\operatorname*{sub}%
\nolimits_{\left(  1,2,\ldots,k\right)  }^{w\left(  Q\right)  }A\text{ and
}a_{i,q_{j}}\\
\text{instead of }n\text{, }A\text{ and }a_{i,j}%
\end{array}
\right) \nonumber\\
&  =\sum_{\alpha\in S_{k}}\left(  -1\right)  ^{\alpha}\prod_{i=1}%
^{k}a_{i,q_{\alpha\left(  i\right)  }}
\label{pf.lem.sol.det.laplace-multi.1.det1}%
\end{align}
(here, we have renamed the summation index $\sigma$ as $\alpha$).

Also, $w\left(  \widetilde{Q}\right)  =\left(  r_{1},r_{2},\ldots
,r_{n-k}\right)  $. Thus,%
\begin{align*}
\operatorname*{sub}\nolimits_{\left(  k+1,k+2,\ldots,n\right)  }^{w\left(
\widetilde{Q}\right)  }B  &  =\operatorname*{sub}\nolimits_{\left(
k+1,k+2,\ldots,n\right)  }^{\left(  r_{1},r_{2},\ldots,r_{n-k}\right)
}B=\operatorname*{sub}\nolimits_{k+1,k+2,\ldots,n}^{r_{1},r_{2},\ldots
,r_{n-k}}B=\left(  b_{k+x,r_{y}}\right)  _{1\leq x\leq n-k,\ 1\leq y\leq
n-k}\\
&  \ \ \ \ \ \ \ \ \ \ \left(
\begin{array}
[c]{c}%
\text{by the definition of }\operatorname*{sub}\nolimits_{k+1,k+2,\ldots
,n}^{r_{1},r_{2},\ldots,r_{n-k}}B\text{,}\\
\text{since }B=\left(  b_{i,j}\right)  _{1\leq i\leq n,\ 1\leq j\leq n}%
\end{array}
\right) \\
&  =\left(  b_{k+i,r_{j}}\right)  _{1\leq i\leq n-k,\ 1\leq j\leq n-k}%
\end{align*}
(here, we have renamed the index $\left(  x,y\right)  $ as $\left(
i,j\right)  $). Thus,%
\begin{align}
\det\left(  \operatorname*{sub}\nolimits_{\left(  k+1,k+2,\ldots,n\right)
}^{w\left(  \widetilde{Q}\right)  }B\right)   &  =\sum_{\sigma\in S_{n-k}%
}\left(  -1\right)  ^{\sigma}\prod_{i=1}^{n-k}b_{k+i,r_{\sigma\left(
i\right)  }}\nonumber\\
&  \ \ \ \ \ \ \ \ \ \ \left(
\begin{array}
[c]{c}%
\text{by (\ref{eq.det.eq.2}), applied to }n-k\text{, }\operatorname*{sub}%
\nolimits_{\left(  k+1,k+2,\ldots,n\right)  }^{w\left(  \widetilde{Q}\right)
}B\text{ and }b_{k+i,r_{j}}\\
\text{instead of }n\text{, }A\text{ and }a_{i,j}%
\end{array}
\right) \nonumber\\
&  =\sum_{\beta\in S_{n-k}}\left(  -1\right)  ^{\beta}\prod_{i=1}%
^{n-k}b_{k+i,r_{\beta\left(  i\right)  }}
\label{pf.lem.sol.det.laplace-multi.1.det2}%
\end{align}
(here, we have renamed the summation index $\sigma$ as $\beta$).

\begin{vershort}
Now, we claim the following: For any $\alpha\in S_{k}$ and $\beta\in S_{n-k}$,
we have%
\begin{equation}
\prod_{i\in\left\{  1,2,\ldots,k\right\}  }a_{i,\sigma_{Q,\alpha,\beta}\left(
i\right)  }=\prod_{i=1}^{k}a_{i,q_{\alpha\left(  i\right)  }}
\label{pf.lem.sol.det.laplace-multi.1.short.factor1}%
\end{equation}
and%
\begin{equation}
\prod_{i\in\left\{  k+1,k+2,\ldots,n\right\}  }b_{i,\sigma_{Q,\alpha,\beta
}\left(  i\right)  }=\prod_{i=1}^{n-k}b_{k+i,r_{\beta\left(  i\right)  }}.
\label{pf.lem.sol.det.laplace-multi.1.short.factor2}%
\end{equation}

[\textit{Proof of (\ref{pf.lem.sol.det.laplace-multi.1.short.factor1}) and
(\ref{pf.lem.sol.det.laplace-multi.1.short.factor2}):} Let $\alpha\in S_{k}$
and $\beta\in S_{n-k}$. The permutation $\sigma_{Q,\alpha,\beta}$ was defined
as the unique $\sigma\in S_{n}$ satisfying%
\begin{equation}
\left(  \sigma\left(  1\right)  ,\sigma\left(  2\right)  ,\ldots,\sigma\left(
n\right)  \right)  =\left(  q_{\alpha\left(  1\right)  },q_{\alpha\left(
2\right)  },\ldots,q_{\alpha\left(  k\right)  },r_{\beta\left(  1\right)
},r_{\beta\left(  2\right)  },\ldots,r_{\beta\left(  n-k\right)  }\right)  .
\label{pf.lem.sol.det.laplace-multi.1.short.factor1.pf.1}%
\end{equation}
Hence, $\sigma_{Q,\alpha,\beta}$ is a $\sigma\in S_{n}$ satisfying
(\ref{pf.lem.sol.det.laplace-multi.1.short.factor1.pf.1}). In other words,
$\sigma_{Q,\alpha,\beta}$ is an element of $S_{n}$ and satisfies%
\[
\left(  \sigma_{Q,\alpha,\beta}\left(  1\right)  ,\sigma_{Q,\alpha,\beta
}\left(  2\right)  ,\ldots,\sigma_{Q,\alpha,\beta}\left(  n\right)  \right)
=\left(  q_{\alpha\left(  1\right)  },q_{\alpha\left(  2\right)  }%
,\ldots,q_{\alpha\left(  k\right)  },r_{\beta\left(  1\right)  }%
,r_{\beta\left(  2\right)  },\ldots,r_{\beta\left(  n-k\right)  }\right)  .
\]
In other words,%
\begin{equation}
\left(  \sigma_{Q,\alpha,\beta}\left(  i\right)  =q_{\alpha\left(  i\right)
}\ \ \ \ \ \ \ \ \ \ \text{for every }i\in\left\{  1,2,\ldots,k\right\}
\right)  \label{pf.lem.sol.det.laplace-multi.1.short.factor1.pf.4}%
\end{equation}
and%
\begin{equation}
\left(  \sigma_{Q,\alpha,\beta}\left(  i\right)  =r_{\beta\left(  i-k\right)
}\ \ \ \ \ \ \ \ \ \ \text{for every }i\in\left\{  k+1,k+2,\ldots,n\right\}
\right)  . \label{pf.lem.sol.det.laplace-multi.1.short.factor2.pf.4}%
\end{equation}
Now,%
\[
\underbrace{\prod_{i\in\left\{  1,2,\ldots,k\right\}  }}_{=\prod_{i=1}^{k}%
}\underbrace{a_{i,\sigma_{Q,\alpha,\beta}\left(  i\right)  }}%
_{\substack{=a_{i,q_{\alpha\left(  i\right)  }}\\\text{(since }\sigma
_{Q,\alpha,\beta}\left(  i\right)  =q_{\alpha\left(  i\right)  }\\\text{(by
(\ref{pf.lem.sol.det.laplace-multi.1.short.factor1.pf.4})))}}}=\prod_{i=1}%
^{k}a_{i,q_{\alpha\left(  i\right)  }}.
\]
This proves (\ref{pf.lem.sol.det.laplace-multi.1.short.factor1}). Furthermore,%
\[
\underbrace{\prod_{i\in\left\{  k+1,k+2,\ldots,n\right\}  }}_{=\prod
_{i=k+1}^{n}}\underbrace{b_{i,\sigma_{Q,\alpha,\beta}\left(  i\right)  }%
}_{\substack{=b_{i,r_{\beta\left(  i-k\right)  }}\\\text{(since }%
\sigma_{Q,\alpha,\beta}\left(  i\right)  =r_{\beta\left(  i-k\right)
}\\\text{(by (\ref{pf.lem.sol.det.laplace-multi.1.short.factor2.pf.4})))}%
}}=\prod_{i=k+1}^{n}b_{i,r_{\beta\left(  i-k\right)  }}=\prod_{i=1}%
^{n-k}b_{k+i,r_{\beta\left(  i\right)  }}%
\]
(here, we have substituted $k+i$ for $i$ in the product). This proves
(\ref{pf.lem.sol.det.laplace-multi.1.short.factor2}).]

Now,%
\begin{align*}
&  \underbrace{\sum_{\substack{\sigma\in S_{n};\\\sigma\left(  \left\{
1,2,\ldots,k\right\}  \right)  =Q}}}_{=\sum_{\sigma\in\left\{  \tau\in
S_{n}\ \mid\ \tau\left(  \left\{  1,2,\ldots,k\right\}  \right)  =Q\right\}
}}\left(  -1\right)  ^{\sigma}\left(  \prod_{i\in\left\{  1,2,\ldots
,k\right\}  }a_{i,\sigma\left(  i\right)  }\right)  \left(  \prod
_{i\in\left\{  k+1,k+2,\ldots,n\right\}  }b_{i,\sigma\left(  i\right)
}\right) \\
&  =\sum_{\sigma\in\left\{  \tau\in S_{n}\ \mid\ \tau\left(  \left\{
1,2,\ldots,k\right\}  \right)  =Q\right\}  }\left(  -1\right)  ^{\sigma
}\left(  \prod_{i\in\left\{  1,2,\ldots,k\right\}  }a_{i,\sigma\left(
i\right)  }\right)  \left(  \prod_{i\in\left\{  k+1,k+2,\ldots,n\right\}
}b_{i,\sigma\left(  i\right)  }\right) \\
&  =\underbrace{\sum_{\left(  \alpha,\beta\right)  \in S_{k}\times S_{n-k}}%
}_{=\sum_{\alpha\in S_{k}}\sum_{\beta\in S_{n-k}}}\left(  -1\right)
^{\sigma_{Q,\alpha,\beta}}\underbrace{\left(  \prod_{i\in\left\{
1,2,\ldots,k\right\}  }a_{i,\sigma_{Q,\alpha,\beta}\left(  i\right)  }\right)
}_{\substack{=\prod_{i=1}^{k}a_{i,q_{\alpha\left(  i\right)  }}\\\text{(by
(\ref{pf.lem.sol.det.laplace-multi.1.short.factor1}))}}}\underbrace{\left(
\prod_{i\in\left\{  k+1,k+2,\ldots,n\right\}  }b_{i,\sigma_{Q,\alpha,\beta
}\left(  i\right)  }\right)  }_{\substack{=\prod_{i=1}^{n-k}b_{k+i,r_{\beta
\left(  i\right)  }}\\\text{(by
(\ref{pf.lem.sol.det.laplace-multi.1.short.factor2}))}}}\\
&  \ \ \ \ \ \ \ \ \ \ \left(
\begin{array}
[c]{c}%
\text{here, we have substituted }\sigma_{Q,\alpha,\beta}\text{ for }%
\sigma\text{ in the sum, since the}\\
\text{map }S_{k}\times S_{n-k}\rightarrow\left\{  \tau\in S_{n}\ \mid
\ \tau\left(  \left\{  1,2,\ldots,k\right\}  \right)  =Q\right\}  ,\ \left(
\alpha,\beta\right)  \mapsto\sigma_{Q,\alpha,\beta}\\
\text{is a bijection}%
\end{array}
\right) \\
&  =\sum_{\alpha\in S_{k}}\sum_{\beta\in S_{n-k}}\underbrace{\left(
-1\right)  ^{\sigma_{Q,\alpha,\beta}}}_{\substack{=\left(  -1\right)
^{\alpha}\cdot\left(  -1\right)  ^{\beta}\cdot\left(  -1\right)  ^{\sum
Q-\left(  1+2+\cdots+k\right)  }\\\text{(by
(\ref{pf.lem.sol.det.laplace-multi.1.sign}))}}}\left(  \prod_{i=1}%
^{k}a_{i,q_{\alpha\left(  i\right)  }}\right)  \left(  \prod_{i=1}%
^{n-k}b_{k+i,r_{\beta\left(  i\right)  }}\right) \\
&  =\sum_{\alpha\in S_{k}}\sum_{\beta\in S_{n-k}}\left(  -1\right)  ^{\alpha
}\cdot\left(  -1\right)  ^{\beta}\cdot\left(  -1\right)  ^{\sum Q-\left(
1+2+\cdots+k\right)  }\left(  \prod_{i=1}^{k}a_{i,q_{\alpha\left(  i\right)
}}\right)  \left(  \prod_{i=1}^{n-k}b_{k+i,r_{\beta\left(  i\right)  }%
}\right)
\end{align*}%
\begin{align*}
&  =\underbrace{\left(  -1\right)  ^{\sum Q-\left(  1+2+\cdots+k\right)  }%
}_{\substack{=\left(  -1\right)  ^{\left(  1+2+\cdots+k\right)  +\sum
Q}\\\text{(since }\sum Q-\left(  1+2+\cdots+k\right)  \\\equiv\sum Q+\left(
1+2+\cdots+k\right)  \\=\left(  1+2+\cdots+k\right)  +\sum Q\operatorname{mod}%
2\text{)}}}\underbrace{\left(  \sum_{\alpha\in S_{k}}\left(  -1\right)
^{\alpha}\prod_{i=1}^{k}a_{i,q_{\alpha\left(  i\right)  }}\right)
}_{\substack{=\det\left(  \operatorname*{sub}\nolimits_{\left(  1,2,\ldots
,k\right)  }^{w\left(  Q\right)  }A\right)  \\\text{(by
(\ref{pf.lem.sol.det.laplace-multi.1.det1}))}}}\underbrace{\left(  \sum
_{\beta\in S_{n-k}}\left(  -1\right)  ^{\beta}\prod_{i=1}^{n-k}b_{k+i,r_{\beta
\left(  i\right)  }}\right)  }_{\substack{=\det\left(  \operatorname*{sub}%
\nolimits_{\left(  k+1,k+2,\ldots,n\right)  }^{w\left(  \widetilde{Q}\right)
}B\right)  \\\text{(by (\ref{pf.lem.sol.det.laplace-multi.1.det2}))}}}\\
&  =\left(  -1\right)  ^{\left(  1+2+\cdots+k\right)  +\sum Q}\det\left(
\operatorname*{sub}\nolimits_{\left(  1,2,\ldots,k\right)  }^{w\left(
Q\right)  }A\right)  \det\left(  \operatorname*{sub}\nolimits_{\left(
k+1,k+2,\ldots,n\right)  }^{w\left(  \widetilde{Q}\right)  }B\right)  .
\end{align*}
This proves Lemma \ref{lem.sol.det.laplace-multi.1}. \qedhere

\end{vershort}

\begin{verlong}
For any $\alpha\in S_{k}$ and $\beta\in S_{n-k}$, we have%
\begin{equation}
\prod_{i\in\left\{  1,2,\ldots,k\right\}  }a_{i,\sigma_{Q,\alpha,\beta}\left(
i\right)  }=\prod_{i=1}^{k}a_{i,q_{\alpha\left(  i\right)  }}
\label{pf.lem.sol.det.laplace-multi.1.factor1}%
\end{equation}
\footnote{\textit{Proof of (\ref{pf.lem.sol.det.laplace-multi.1.factor1}):}
Let $\alpha\in S_{k}$ and $\beta\in S_{n-k}$. Now, $\sigma_{Q,\alpha,\beta}$
is the unique $\sigma\in S_{n}$ satisfying%
\begin{equation}
\left(  \sigma\left(  1\right)  ,\sigma\left(  2\right)  ,\ldots,\sigma\left(
n\right)  \right)  =\left(  q_{\alpha\left(  1\right)  },q_{\alpha\left(
2\right)  },\ldots,q_{\alpha\left(  k\right)  },r_{\beta\left(  1\right)
},r_{\beta\left(  2\right)  },\ldots,r_{\beta\left(  n-k\right)  }\right)
\label{pf.lem.sol.det.laplace-multi.1.factor1.pf.1}%
\end{equation}
(because this is how $\sigma_{Q,\alpha,\beta}$ is defined). Hence,
$\sigma_{Q,\alpha,\beta}$ is a $\sigma\in S_{n}$ satisfying
(\ref{pf.lem.sol.det.laplace-multi.1.factor1.pf.1}). In other words,
$\sigma_{Q,\alpha,\beta}$ is an element of $S_{n}$ and satisfies%
\begin{align}
&  \left(  \sigma_{Q,\alpha,\beta}\left(  1\right)  ,\sigma_{Q,\alpha,\beta
}\left(  2\right)  ,\ldots,\sigma_{Q,\alpha,\beta}\left(  n\right)  \right)
\nonumber\\
&  =\left(  q_{\alpha\left(  1\right)  },q_{\alpha\left(  2\right)  }%
,\ldots,q_{\alpha\left(  k\right)  },r_{\beta\left(  1\right)  }%
,r_{\beta\left(  2\right)  },\ldots,r_{\beta\left(  n-k\right)  }\right)  .
\label{pf.lem.sol.det.laplace-multi.1.factor1.pf.2}%
\end{align}
\par
Now,%
\begin{align*}
&  \left(  \sigma_{Q,\alpha,\beta}\left(  1\right)  ,\sigma_{Q,\alpha,\beta
}\left(  2\right)  ,\ldots,\sigma_{Q,\alpha,\beta}\left(  k\right)  \right) \\
&  =\left(  \text{the list of the first }k\text{ entries of the list
}\underbrace{\left(  \sigma_{Q,\alpha,\beta}\left(  1\right)  ,\sigma
_{Q,\alpha,\beta}\left(  2\right)  ,\ldots,\sigma_{Q,\alpha,\beta}\left(
n\right)  \right)  }_{\substack{=\left(  q_{\alpha\left(  1\right)
},q_{\alpha\left(  2\right)  },\ldots,q_{\alpha\left(  k\right)  }%
,r_{\beta\left(  1\right)  },r_{\beta\left(  2\right)  },\ldots,r_{\beta
\left(  n-k\right)  }\right)  \\\text{(by
(\ref{pf.lem.sol.det.laplace-multi.1.factor1.pf.2}))}}}\right) \\
&  =\left(  \text{the list of the first }k\text{ entries of the list }\left(
q_{\alpha\left(  1\right)  },q_{\alpha\left(  2\right)  },\ldots
,q_{\alpha\left(  k\right)  },r_{\beta\left(  1\right)  },r_{\beta\left(
2\right)  },\ldots,r_{\beta\left(  n-k\right)  }\right)  \right) \\
&  =\left(  q_{\alpha\left(  1\right)  },q_{\alpha\left(  2\right)  }%
,\ldots,q_{\alpha\left(  k\right)  }\right)  .
\end{align*}
In other words,%
\begin{equation}
\sigma_{Q,\alpha,\beta}\left(  i\right)  =q_{\alpha\left(  i\right)
}\ \ \ \ \ \ \ \ \ \ \text{for every }i\in\left\{  1,2,\ldots,k\right\}  .
\label{pf.lem.sol.det.laplace-multi.1.factor1.pf.4}%
\end{equation}
Hence,%
\[
\underbrace{\prod_{i\in\left\{  1,2,\ldots,k\right\}  }}_{=\prod_{i=1}^{k}%
}\underbrace{a_{i,\sigma_{Q,\alpha,\beta}\left(  i\right)  }}%
_{\substack{=a_{i,q_{\alpha\left(  i\right)  }}\\\text{(since }\sigma
_{Q,\alpha,\beta}\left(  i\right)  =q_{\alpha\left(  i\right)  }\\\text{(by
(\ref{pf.lem.sol.det.laplace-multi.1.factor1.pf.4})))}}}=\prod_{i=1}%
^{k}a_{i,q_{\alpha\left(  i\right)  }}.
\]
This proves (\ref{pf.lem.sol.det.laplace-multi.1.factor1}).} and%
\begin{equation}
\prod_{i\in\left\{  k+1,k+2,\ldots,n\right\}  }b_{i,\sigma_{Q,\alpha,\beta
}\left(  i\right)  }=\prod_{i=1}^{n-k}b_{k+i,r_{\beta\left(  i\right)  }}
\label{pf.lem.sol.det.laplace-multi.1.factor2}%
\end{equation}
\footnote{\textit{Proof of (\ref{pf.lem.sol.det.laplace-multi.1.factor2}):}
Let $\alpha\in S_{k}$ and $\beta\in S_{n-k}$. Now, $\sigma_{Q,\alpha,\beta}$
is the unique $\sigma\in S_{n}$ satisfying%
\begin{equation}
\left(  \sigma\left(  1\right)  ,\sigma\left(  2\right)  ,\ldots,\sigma\left(
n\right)  \right)  =\left(  q_{\alpha\left(  1\right)  },q_{\alpha\left(
2\right)  },\ldots,q_{\alpha\left(  k\right)  },r_{\beta\left(  1\right)
},r_{\beta\left(  2\right)  },\ldots,r_{\beta\left(  n-k\right)  }\right)
\label{pf.lem.sol.det.laplace-multi.1.factor2.pf.1}%
\end{equation}
(because this is how $\sigma_{Q,\alpha,\beta}$ is defined). Hence,
$\sigma_{Q,\alpha,\beta}$ is a $\sigma\in S_{n}$ satisfying
(\ref{pf.lem.sol.det.laplace-multi.1.factor2.pf.1}). In other words,
$\sigma_{Q,\alpha,\beta}$ is an element of $S_{n}$ and satisfies%
\begin{align}
&  \left(  \sigma_{Q,\alpha,\beta}\left(  1\right)  ,\sigma_{Q,\alpha,\beta
}\left(  2\right)  ,\ldots,\sigma_{Q,\alpha,\beta}\left(  n\right)  \right)
\nonumber\\
&  =\left(  q_{\alpha\left(  1\right)  },q_{\alpha\left(  2\right)  }%
,\ldots,q_{\alpha\left(  k\right)  },r_{\beta\left(  1\right)  }%
,r_{\beta\left(  2\right)  },\ldots,r_{\beta\left(  n-k\right)  }\right)  .
\label{pf.lem.sol.det.laplace-multi.1.factor2.pf.2}%
\end{align}
\par
Now,%
\begin{align*}
&  \left(  \sigma_{Q,\alpha,\beta}\left(  k+1\right)  ,\sigma_{Q,\alpha,\beta
}\left(  k+2\right)  ,\ldots,\sigma_{Q,\alpha,\beta}\left(  n\right)  \right)
\\
&  =\left(  \text{the list of the last }n-k\text{ entries of the list
}\underbrace{\left(  \sigma_{Q,\alpha,\beta}\left(  1\right)  ,\sigma
_{Q,\alpha,\beta}\left(  2\right)  ,\ldots,\sigma_{Q,\alpha,\beta}\left(
n\right)  \right)  }_{\substack{=\left(  q_{\alpha\left(  1\right)
},q_{\alpha\left(  2\right)  },\ldots,q_{\alpha\left(  k\right)  }%
,r_{\beta\left(  1\right)  },r_{\beta\left(  2\right)  },\ldots,r_{\beta
\left(  n-k\right)  }\right)  \\\text{(by
(\ref{pf.lem.sol.det.laplace-multi.1.factor2.pf.2}))}}}\right) \\
&  =\left(  \text{the list of the last }n-k\text{ entries of the list }\left(
q_{\alpha\left(  1\right)  },q_{\alpha\left(  2\right)  },\ldots
,q_{\alpha\left(  k\right)  },r_{\beta\left(  1\right)  },r_{\beta\left(
2\right)  },\ldots,r_{\beta\left(  n-k\right)  }\right)  \right) \\
&  =\left(  r_{\beta\left(  1\right)  },r_{\beta\left(  2\right)  }%
,\ldots,r_{\beta\left(  n-k\right)  }\right)  .
\end{align*}
In other words,%
\begin{equation}
\sigma_{Q,\alpha,\beta}\left(  k+i\right)  =r_{\beta\left(  i\right)
}\ \ \ \ \ \ \ \ \ \ \text{for every }i\in\left\{  1,2,\ldots,n-k\right\}  .
\label{pf.lem.sol.det.laplace-multi.1.factor2.pf.4}%
\end{equation}
Now,%
\begin{align*}
\underbrace{\prod_{i\in\left\{  k+1,k+2,\ldots,n\right\}  }}_{=\prod
_{i=k+1}^{n}}b_{i,\sigma_{Q,\alpha,\beta}\left(  i\right)  }  &
=\prod_{i=k+1}^{n}b_{i,\sigma_{Q,\alpha,\beta}\left(  i\right)  }=\prod
_{i=1}^{n-k}\underbrace{b_{i+k,\sigma_{Q,\alpha,\beta}\left(  k+i\right)  }%
}_{\substack{=b_{k+i,\sigma_{Q,\alpha,\beta}\left(  k+i\right)  }%
\\\text{(since }i+k=k+i\text{)}}}\\
&  \ \ \ \ \ \ \ \ \ \ \left(  \text{here, we have substituted }i+k\text{ for
}i\text{ in the product}\right) \\
&  =\prod_{i=1}^{n-k}\underbrace{b_{k+i,\sigma_{Q,\alpha,\beta}\left(
k+i\right)  }}_{\substack{=b_{k+i,r_{\beta\left(  i\right)  }}\\\text{(since
}\sigma_{Q,\alpha,\beta}\left(  k+i\right)  =r_{\beta\left(  i\right)
}\\\text{(by (\ref{pf.lem.sol.det.laplace-multi.1.factor2.pf.4})))}}%
}=\prod_{i=1}^{n-k}b_{k+i,r_{\beta\left(  i\right)  }}.
\end{align*}
This proves (\ref{pf.lem.sol.det.laplace-multi.1.factor2}).}.

Now,%
\begin{align*}
&  \underbrace{\sum_{\substack{\sigma\in S_{n};\\\sigma\left(  \left\{
1,2,\ldots,k\right\}  \right)  =Q}}}_{=\sum_{\sigma\in\left\{  \tau\in
S_{n}\ \mid\ \tau\left(  \left\{  1,2,\ldots,k\right\}  \right)  =Q\right\}
}}\left(  -1\right)  ^{\sigma}\left(  \prod_{i\in\left\{  1,2,\ldots
,k\right\}  }a_{i,\sigma\left(  i\right)  }\right)  \left(  \prod
_{i\in\left\{  k+1,k+2,\ldots,n\right\}  }b_{i,\sigma\left(  i\right)
}\right) \\
&  =\sum_{\sigma\in\left\{  \tau\in S_{n}\ \mid\ \tau\left(  \left\{
1,2,\ldots,k\right\}  \right)  =Q\right\}  }\left(  -1\right)  ^{\sigma
}\left(  \prod_{i\in\left\{  1,2,\ldots,k\right\}  }a_{i,\sigma\left(
i\right)  }\right)  \left(  \prod_{i\in\left\{  k+1,k+2,\ldots,n\right\}
}b_{i,\sigma\left(  i\right)  }\right) \\
&  =\underbrace{\sum_{\left(  \alpha,\beta\right)  \in S_{k}\times S_{n-k}}%
}_{=\sum_{\alpha\in S_{k}}\sum_{\beta\in S_{n-k}}}\left(  -1\right)
^{\sigma_{Q,\alpha,\beta}}\underbrace{\left(  \prod_{i\in\left\{
1,2,\ldots,k\right\}  }a_{i,\sigma_{Q,\alpha,\beta}\left(  i\right)  }\right)
}_{\substack{=\prod_{i=1}^{k}a_{i,q_{\alpha\left(  i\right)  }}\\\text{(by
(\ref{pf.lem.sol.det.laplace-multi.1.factor1}))}}}\underbrace{\left(
\prod_{i\in\left\{  k+1,k+2,\ldots,n\right\}  }b_{i,\sigma_{Q,\alpha,\beta
}\left(  i\right)  }\right)  }_{\substack{=\prod_{i=1}^{n-k}b_{k+i,r_{\beta
\left(  i\right)  }}\\\text{(by (\ref{pf.lem.sol.det.laplace-multi.1.factor2}%
))}}}\\
&  \ \ \ \ \ \ \ \ \ \ \left(
\begin{array}
[c]{c}%
\text{here, we have substituted }\sigma_{Q,\alpha,\beta}\text{ for }%
\sigma\text{ in the sum, since the}\\
\text{map }S_{k}\times S_{n-k}\rightarrow\left\{  \tau\in S_{n}\ \mid
\ \tau\left(  \left\{  1,2,\ldots,k\right\}  \right)  =Q\right\}  ,\ \left(
\alpha,\beta\right)  \mapsto\sigma_{Q,\alpha,\beta}\\
\text{is a bijection}%
\end{array}
\right) \\
&  =\sum_{\alpha\in S_{k}}\sum_{\beta\in S_{n-k}}\underbrace{\left(
-1\right)  ^{\sigma_{Q,\alpha,\beta}}}_{\substack{=\left(  -1\right)
^{\alpha}\cdot\left(  -1\right)  ^{\beta}\cdot\left(  -1\right)  ^{\sum
Q-\left(  1+2+\cdots+k\right)  }\\\text{(by
(\ref{pf.lem.sol.det.laplace-multi.1.sign}))}}}\left(  \prod_{i=1}%
^{k}a_{i,q_{\alpha\left(  i\right)  }}\right)  \left(  \prod_{i=1}%
^{n-k}b_{k+i,r_{\beta\left(  i\right)  }}\right) \\
&  =\sum_{\alpha\in S_{k}}\sum_{\beta\in S_{n-k}}\left(  -1\right)  ^{\alpha
}\cdot\left(  -1\right)  ^{\beta}\cdot\left(  -1\right)  ^{\sum Q-\left(
1+2+\cdots+k\right)  }\left(  \prod_{i=1}^{k}a_{i,q_{\alpha\left(  i\right)
}}\right)  \left(  \prod_{i=1}^{n-k}b_{k+i,r_{\beta\left(  i\right)  }%
}\right)
\end{align*}%
\begin{align*}
&  =\left(  -1\right)  ^{\sum Q-\left(  1+2+\cdots+k\right)  }\underbrace{\sum
_{\alpha\in S_{k}}\sum_{\beta\in S_{n-k}}\left(  -1\right)  ^{\alpha}%
\cdot\left(  -1\right)  ^{\beta}\cdot\left(  \prod_{i=1}^{k}a_{i,q_{\alpha
\left(  i\right)  }}\right)  \left(  \prod_{i=1}^{n-k}b_{k+i,r_{\beta\left(
i\right)  }}\right)  }_{=\left(  \sum_{\alpha\in S_{k}}\left(  -1\right)
^{\alpha}\prod_{i=1}^{k}a_{i,q_{\alpha\left(  i\right)  }}\right)  \left(
\sum_{\beta\in S_{n-k}}\left(  -1\right)  ^{\beta}\prod_{i=1}^{n-k}%
b_{k+i,r_{\beta\left(  i\right)  }}\right)  }\\
&  =\underbrace{\left(  -1\right)  ^{\sum Q-\left(  1+2+\cdots+k\right)  }%
}_{\substack{=\left(  -1\right)  ^{\left(  1+2+\cdots+k\right)  +\sum
Q}\\\text{(since }\sum Q-\left(  1+2+\cdots+k\right)  \\\equiv\sum Q+\left(
1+2+\cdots+k\right)  \\=\left(  1+2+\cdots+k\right)  +\sum Q\operatorname{mod}%
2\text{)}}}\underbrace{\left(  \sum_{\alpha\in S_{k}}\left(  -1\right)
^{\alpha}\prod_{i=1}^{k}a_{i,q_{\alpha\left(  i\right)  }}\right)
}_{\substack{=\det\left(  \operatorname*{sub}\nolimits_{\left(  1,2,\ldots
,k\right)  }^{w\left(  Q\right)  }A\right)  \\\text{(by
(\ref{pf.lem.sol.det.laplace-multi.1.det1}))}}}\underbrace{\left(  \sum
_{\beta\in S_{n-k}}\left(  -1\right)  ^{\beta}\prod_{i=1}^{n-k}b_{k+i,r_{\beta
\left(  i\right)  }}\right)  }_{\substack{=\det\left(  \operatorname*{sub}%
\nolimits_{\left(  k+1,k+2,\ldots,n\right)  }^{w\left(  \widetilde{Q}\right)
}B\right)  \\\text{(by (\ref{pf.lem.sol.det.laplace-multi.1.det2}))}}}\\
&  =\left(  -1\right)  ^{\left(  1+2+\cdots+k\right)  +\sum Q}\det\left(
\operatorname*{sub}\nolimits_{\left(  1,2,\ldots,k\right)  }^{w\left(
Q\right)  }A\right)  \det\left(  \operatorname*{sub}\nolimits_{\left(
k+1,k+2,\ldots,n\right)  }^{w\left(  \widetilde{Q}\right)  }B\right)  .
\end{align*}
This proves Lemma \ref{lem.sol.det.laplace-multi.1}.
\end{verlong}
\end{proof}

Before we move on to the proof of Lemma \ref{lem.det.laplace-multi.Apq}, let
us show two more lemmas. The first one is a really simple fact about symmetric groups:

\begin{lemma}
\label{lem.sol.det.laplace-multi.group-bij}Let $n\in\mathbb{N}$. Let
$\gamma\in S_{n}$. Then, the map $S_{n}\rightarrow S_{n},\ \sigma\mapsto
\sigma\circ\gamma$ is a bijection.
\end{lemma}

\begin{vershort}
\begin{proof}
[Proof of Lemma \ref{lem.sol.det.laplace-multi.group-bij}.]It is easy to see
that the maps $S_{n}\rightarrow S_{n},\ \sigma\mapsto\sigma\circ\gamma$ and
$S_{n}\rightarrow S_{n},\ \sigma\mapsto\sigma\circ\gamma^{-1}$ are mutually
inverse. Thus, the map $S_{n}\rightarrow S_{n},\ \sigma\mapsto\sigma
\circ\gamma$ is invertible, i.e., a bijection. Lemma
\ref{lem.sol.det.laplace-multi.group-bij} is proven.
\end{proof}
\end{vershort}

\begin{verlong}
\begin{proof}
[Proof of Lemma \ref{lem.sol.det.laplace-multi.group-bij}.]Let $\mathbf{A}$ be
the map $S_{n}\rightarrow S_{n},\ \sigma\mapsto\sigma\circ\gamma$. Let
$\mathbf{B}$ be the map $S_{n}\rightarrow S_{n},\ \sigma\mapsto\sigma
\circ\gamma^{-1}$. Then, $\mathbf{A}\circ\mathbf{B}=\operatorname*{id}%
\nolimits_{S_{n}}$\ \ \ \ \footnote{\textit{Proof.} Every $\sigma\in S_{n}$
satisfies%
\begin{align*}
\left(  \mathbf{A}\circ\mathbf{B}\right)  \left(  \sigma\right)   &
=\mathbf{A}\left(  \underbrace{\mathbf{B}\left(  \sigma\right)  }%
_{\substack{=\sigma\circ\gamma^{-1}\\\text{(by the definition of }%
\mathbf{B}\text{)}}}\right)  =\mathbf{A}\left(  \sigma\circ\gamma^{-1}\right)
=\left(  \sigma\circ\gamma^{-1}\right)  \circ\gamma\ \ \ \ \ \ \ \ \ \ \left(
\text{by the definition of }\mathbf{A}\right) \\
&  =\sigma\circ\underbrace{\gamma^{-1}\circ\gamma}_{=\operatorname*{id}%
}=\sigma\circ\operatorname*{id}=\sigma=\operatorname*{id}\nolimits_{S_{n}%
}\left(  \sigma\right)  .
\end{align*}
In other words, $\mathbf{A}\circ\mathbf{B}=\operatorname*{id}\nolimits_{S_{n}%
}$. Qed.} and $\mathbf{B}\circ\mathbf{A}=\operatorname*{id}\nolimits_{S_{n}}%
$\ \ \ \ \footnote{\textit{Proof.} Every $\sigma\in S_{n}$ satisfies%
\begin{align*}
\left(  \mathbf{B}\circ\mathbf{A}\right)  \left(  \sigma\right)   &
=\mathbf{B}\left(  \underbrace{\mathbf{A}\left(  \sigma\right)  }%
_{\substack{=\sigma\circ\gamma\\\text{(by the definition of }\mathbf{A}%
\text{)}}}\right)  =\mathbf{B}\left(  \sigma\circ\gamma\right)  =\left(
\sigma\circ\gamma\right)  \circ\gamma^{-1}\ \ \ \ \ \ \ \ \ \ \left(  \text{by
the definition of }\mathbf{B}\right) \\
&  =\sigma\circ\underbrace{\gamma\circ\gamma^{-1}}_{=\operatorname*{id}%
}=\sigma\circ\operatorname*{id}=\sigma=\operatorname*{id}\nolimits_{S_{n}%
}\left(  \sigma\right)  .
\end{align*}
In other words, $\mathbf{B}\circ\mathbf{A}=\operatorname*{id}\nolimits_{S_{n}%
}$. Qed.}. These two equalities show that the maps $\mathbf{A}$ and
$\mathbf{B}$ are mutually inverse. Hence, the map $\mathbf{A}$ is invertible.
In other words, the map $\mathbf{A}$ is a bijection. In other words, the map
$S_{n}\rightarrow S_{n},\ \sigma\mapsto\sigma\circ\gamma$ is a bijection
(since the map $\mathbf{A}$ is the map $S_{n}\rightarrow S_{n},\ \sigma
\mapsto\sigma\circ\gamma$). This proves Lemma
\ref{lem.sol.det.laplace-multi.group-bij}.
\end{proof}
\end{verlong}

Next, we show a lemma which is a distillate of some parts of Exercise
\ref{exe.Ialbe}:

\begin{lemma}
\label{lem.sol.addexe.jacobi-complement.Ialbe}Let $n\in\mathbb{N}$. For any
subset $I$ of $\left\{  1,2,\ldots,n\right\}  $, we let $\widetilde{I}$ denote
the complement $\left\{  1,2,\ldots,n\right\}  \setminus I$ of $I$.

Let $I$ be a subset of $\left\{  1,2,\ldots,n\right\}  $. Let $k=\left\vert
I\right\vert $. Then, there exists a $\sigma\in S_{n}$ satisfying $\left(
\sigma\left(  1\right)  ,\sigma\left(  2\right)  ,\ldots,\sigma\left(
k\right)  \right)  =w\left(  I\right)  $, $\left(  \sigma\left(  k+1\right)
,\sigma\left(  k+2\right)  ,\ldots,\sigma\left(  n\right)  \right)  =w\left(
\widetilde{I}\right)  $ and $\left(  -1\right)  ^{\sigma}=\left(  -1\right)
^{\sum I-\left(  1+2+\cdots+k\right)  }$.
\end{lemma}

\begin{proof}
[Proof of Lemma \ref{lem.sol.addexe.jacobi-complement.Ialbe}.]We begin by
introducing some notations:

\begin{verlong}
For every $h\in\mathbb{N}$, we let $\operatorname*{id}\nolimits_{h}$ denote
the identity permutation $\operatorname*{id}\nolimits_{\left\{  1,2,\ldots
,h\right\}  }\in S_{h}$. For every $h\in\mathbb{N}$, we have%
\begin{equation}
\left(  -1\right)  ^{\operatorname*{id}\nolimits_{h}}=1.
\label{pf.lem.sol.addexe.jacobi-complement.Ialbe.1}%
\end{equation}
(Indeed, this is just a restatement of the well-known fact that $\left(
-1\right)  ^{\operatorname*{id}}=1$, where $\operatorname*{id}%
=\operatorname*{id}\nolimits_{h}$ is the identity permutation in $S_{h}$.)
\end{verlong}

\begin{vershort}
The definition of $\widetilde{I}$ yields $\widetilde{I}=\left\{
1,2,\ldots,n\right\}  \setminus I$. Since $I$ is a subset of $\left\{
1,2,\ldots,n\right\}  $, this leads to $\left\vert \widetilde{I}\right\vert
=\underbrace{\left\vert \left\{  1,2,\ldots,n\right\}  \right\vert }%
_{=n}-\underbrace{\left\vert I\right\vert }_{=k}=n-k$. Thus, $n-k=\left\vert
\widetilde{I}\right\vert \geq0$, so that $n\geq k$ and thus $k\in\left\{
0,1,\ldots,n\right\}  $.
\end{vershort}

\begin{verlong}
The definition of $\widetilde{I}$ yields $\widetilde{I}=\left\{
1,2,\ldots,n\right\}  \setminus I\subseteq\left\{  1,2,\ldots,n\right\}  $.
Thus, $\widetilde{I}$ is a finite set (since $\left\{  1,2,\ldots,n\right\}  $
is a finite set). Hence, $\left\vert \widetilde{I}\right\vert \in\mathbb{N}$.

Also, $I\subseteq\left\{  1,2,\ldots,n\right\}  $. Thus, $I$ is a finite set
(since $\left\{  1,2,\ldots,n\right\}  $ is a finite set). Hence, $\left\vert
I\right\vert \in\mathbb{N}$. Thus, $k=\left\vert I\right\vert \in\mathbb{N}$.

Also,%
\begin{align*}
\left\vert \underbrace{\widetilde{I}}_{=\left\{  1,2,\ldots,n\right\}
\setminus I}\right\vert  &  =\left\vert \left\{  1,2,\ldots,n\right\}
\setminus I\right\vert =\underbrace{\left\vert \left\{  1,2,\ldots,n\right\}
\right\vert }_{=n}-\underbrace{\left\vert I\right\vert }_{=k}%
\ \ \ \ \ \ \ \ \ \ \left(  \text{since }I\subseteq\left\{  1,2,\ldots
,n\right\}  \right) \\
&  =n-k.
\end{align*}
Hence, $n-k=\left\vert \widetilde{I}\right\vert \in\mathbb{N}$; thus,
$n-k\geq0$. In other words, $k\leq n$. Combined with $k\geq0$ (since
$k=\left\vert I\right\vert \in\mathbb{N}$), this yields $k\in\left\{
0,1,\ldots,n\right\}  $.
\end{verlong}

We know that $w\left(  I\right)  $ is the list of all elements of $I$ in
increasing order (with no repetitions) (by the definition of $w\left(
I\right)  $). Thus, $w\left(  I\right)  $ is a list of $\left\vert
I\right\vert $ elements. In other words, $w\left(  I\right)  $ is a list of
$k$ elements (since $\left\vert I\right\vert =k$).

Write $w\left(  I\right)  $ in the form $w\left(  I\right)  =\left(
a_{1},a_{2},\ldots,a_{k}\right)  $. (This is possible, since $w\left(
I\right)  $ is a list of $k$ elements.)

We know that $w\left(  I\right)  $ is the list of all elements of $I$ in
increasing order (with no repetitions). In other words, $\left(  a_{1}%
,a_{2},\ldots,a_{k}\right)  $ is the list of all elements of $I$ in increasing
order (with no repetitions) (since $w\left(  I\right)  =\left(  a_{1}%
,a_{2},\ldots,a_{k}\right)  $).

We know that $w\left(  \widetilde{I}\right)  $ is the list of all elements of
$\widetilde{I}$ in increasing order (with no repetitions) (by the definition
of $w\left(  \widetilde{I}\right)  $). Thus, $w\left(  \widetilde{I}\right)  $
is a list of $\left\vert \widetilde{I}\right\vert $ elements. In other words,
$w\left(  \widetilde{I}\right)  $ is a list of $n-k$ elements (since
$\left\vert \widetilde{I}\right\vert =n-k$).

Write $w\left(  \widetilde{I}\right)  $ in the form $w\left(  \widetilde{I}%
\right)  =\left(  b_{1},b_{2},\ldots,b_{n-k}\right)  $. (This is possible,
since $w\left(  \widetilde{I}\right)  $ is a list of $n-k$ elements.)

We know that $w\left(  \widetilde{I}\right)  $ is the list of all elements of
$\widetilde{I}$ in increasing order (with no repetitions). In other words,
$\left(  b_{1},b_{2},\ldots,b_{n-k}\right)  $ is the list of all elements of
$\widetilde{I}$ in increasing order (with no repetitions) (since $w\left(
\widetilde{I}\right)  =\left(  b_{1},b_{2},\ldots,b_{n-k}\right)  $). In other
words, $\left(  b_{1},b_{2},\ldots,b_{n-k}\right)  $ is the list of all
elements of $\left\{  1,2,\ldots,n\right\}  \setminus I$ in increasing order
(with no repetitions) (since $\widetilde{I}=\left\{  1,2,\ldots,n\right\}
\setminus I$).

Now we know that $\left(  a_{1},a_{2},\ldots,a_{k}\right)  $ is the list of
all elements of $I$ in increasing order (with no repetitions), and that
$\left(  b_{1},b_{2},\ldots,b_{n-k}\right)  $ is the list of all elements of
$\left\{  1,2,\ldots,n\right\}  \setminus I$ in increasing order (with no
repetitions). Hence, for every $\alpha\in S_{k}$ and $\beta\in S_{n-k}$, we
can define an element $\sigma_{I,\alpha,\beta}\in S_{n}$ according to Exercise
\ref{exe.Ialbe} \textbf{(a)}. Consider this $\sigma_{I,\alpha,\beta}$.
Exercise \ref{exe.Ialbe} \textbf{(b)} shows that for every $\alpha\in S_{k}$
and $\beta\in S_{n-k}$, we have%
\[
\ell\left(  \sigma_{I,\alpha,\beta}\right)  =\ell\left(  \alpha\right)
+\ell\left(  \beta\right)  +\sum I-\left(  1+2+\cdots+k\right)
\]
and%
\begin{equation}
\left(  -1\right)  ^{\sigma_{I,\alpha,\beta}}=\left(  -1\right)  ^{\alpha
}\cdot\left(  -1\right)  ^{\beta}\cdot\left(  -1\right)  ^{\sum I-\left(
1+2+\cdots+k\right)  }.
\label{pf.lem.sol.addexe.jacobi-complement.Ialbe.sign-gen}%
\end{equation}

\begin{vershort}
Now, consider the identity permutations $\operatorname*{id}\in S_{k}$ and
$\operatorname*{id}\in S_{n-k}$. (These are (in general) two different
permutations, although we denote them both by $\operatorname*{id}$.) They give
rise to an element $\sigma_{I,\operatorname*{id},\operatorname*{id}}\in S_{n}$
(obtained by setting $\alpha=\operatorname*{id}\in S_{k}$ and $\beta
=\operatorname*{id}\in S_{n-k}$ in the definition of $\sigma_{I,\alpha,\beta}%
$). Denote this element $\sigma_{I,\operatorname*{id},\operatorname*{id}}$ by
$\gamma$. Thus, $\gamma=\sigma_{I,\operatorname*{id},\operatorname*{id}}$, so
that%
\begin{align}
\left(  -1\right)  ^{\gamma}  &  =\left(  -1\right)  ^{\sigma
_{I,\operatorname*{id},\operatorname*{id}}}=\underbrace{\left(  -1\right)
^{\operatorname*{id}}}_{=1}\cdot\underbrace{\left(  -1\right)
^{\operatorname*{id}}}_{=1}\cdot\left(  -1\right)  ^{\sum I-\left(
1+2+\cdots+k\right)  }\nonumber\\
&  \ \ \ \ \ \ \ \ \ \ \left(  \text{by
(\ref{pf.lem.sol.addexe.jacobi-complement.Ialbe.sign-gen}) (applied to }%
\alpha=\operatorname*{id}\text{ and }\beta=\operatorname*{id}\text{)}\right)
\nonumber\\
&  =\left(  -1\right)  ^{\sum I-\left(  1+2+\cdots+k\right)  }.
\label{pf.lem.sol.addexe.jacobi-complement.Ialbe.short.sign-gamma}%
\end{align}

\end{vershort}

\begin{verlong}
Now, $\operatorname*{id}\nolimits_{k}\in S_{k}$ and $\operatorname*{id}%
\nolimits_{n-k}\in S_{n-k}$. Hence, an element $\sigma_{I,\operatorname*{id}%
\nolimits_{k},\operatorname*{id}\nolimits_{n-k}}\in S_{n}$ is defined (because
we can set $\alpha=\operatorname*{id}\nolimits_{k}$ and $\beta
=\operatorname*{id}\nolimits_{n-k}$ in the definition of $\sigma
_{I,\alpha,\beta}$). Denote this element $\sigma_{I,\operatorname*{id}%
\nolimits_{k},\operatorname*{id}\nolimits_{n-k}}$ by $\gamma$. Thus,
$\gamma=\sigma_{I,\operatorname*{id}\nolimits_{k},\operatorname*{id}%
\nolimits_{n-k}}$, so that%
\begin{align}
\left(  -1\right)  ^{\gamma}  &  =\left(  -1\right)  ^{\sigma
_{I,\operatorname*{id}\nolimits_{k},\operatorname*{id}\nolimits_{n-k}}%
}=\underbrace{\left(  -1\right)  ^{\operatorname*{id}\nolimits_{k}}%
}_{\substack{=1\\\text{(by (\ref{pf.lem.sol.addexe.jacobi-complement.Ialbe.1}%
)}\\\text{(applied to }h=k\text{))}}}\cdot\underbrace{\left(  -1\right)
^{\operatorname*{id}\nolimits_{n-k}}}_{\substack{=1\\\text{(by
(\ref{pf.lem.sol.addexe.jacobi-complement.Ialbe.1})}\\\text{(applied to
}h=n-k\text{))}}}\cdot\left(  -1\right)  ^{\sum I-\left(  1+2+\cdots+k\right)
}\nonumber\\
&  \ \ \ \ \ \ \ \ \ \ \left(  \text{by
(\ref{pf.lem.sol.addexe.jacobi-complement.Ialbe.sign-gen}) (applied to }%
\alpha=\operatorname*{id}\nolimits_{k}\text{ and }\beta=\operatorname*{id}%
\nolimits_{n-k}\text{)}\right) \nonumber\\
&  =\left(  -1\right)  ^{\sum I-\left(  1+2+\cdots+k\right)  }.
\label{pf.lem.sol.addexe.jacobi-complement.Ialbe.sign-gamma}%
\end{align}

\end{verlong}

Now, we claim that%
\begin{equation}
\left(  \gamma\left(  1\right)  ,\gamma\left(  2\right)  ,\ldots,\gamma\left(
k\right)  \right)  =\left(  a_{1},a_{2},\ldots,a_{k}\right)
\label{pf.lem.sol.addexe.jacobi-complement.Ialbe.gamma-tuple1}%
\end{equation}
and%
\begin{equation}
\left(  \gamma\left(  k+1\right)  ,\gamma\left(  k+2\right)  ,\ldots
,\gamma\left(  n\right)  \right)  =\left(  b_{1},b_{2},\ldots,b_{n-k}\right)
. \label{pf.lem.sol.addexe.jacobi-complement.Ialbe.gamma-tuple2}%
\end{equation}

\begin{vershort}
[\textit{Proof of
(\ref{pf.lem.sol.addexe.jacobi-complement.Ialbe.gamma-tuple1}) and
(\ref{pf.lem.sol.addexe.jacobi-complement.Ialbe.gamma-tuple2}):} The
permutation $\gamma$ is $\sigma_{I,\operatorname*{id},\operatorname*{id}}$. In
other words, the permutation $\gamma$ is the unique $\sigma\in S_{n}$
satisfying%
\begin{equation}
\left(  \sigma\left(  1\right)  ,\sigma\left(  2\right)  ,\ldots,\sigma\left(
n\right)  \right)  =\left(  a_{\operatorname*{id}\left(  1\right)
},a_{\operatorname*{id}\left(  2\right)  },\ldots,a_{\operatorname*{id}\left(
k\right)  },b_{\operatorname*{id}\left(  1\right)  },b_{\operatorname*{id}%
\left(  2\right)  },\ldots,b_{\operatorname*{id}\left(  n-k\right)  }\right)
\label{pf.lem.sol.addexe.jacobi-complement.Ialbe.short.gamma-tuple.pf.1}%
\end{equation}
(because this is how $\sigma_{I,\operatorname*{id},\operatorname*{id}}$ was
defined). Hence, $\gamma$ is a $\sigma\in S_{n}$ satisfying
(\ref{pf.lem.sol.addexe.jacobi-complement.Ialbe.short.gamma-tuple.pf.1}). In
other words, $\gamma$ is an element of $S_{n}$ and satisfies%
\[
\left(  \gamma\left(  1\right)  ,\gamma\left(  2\right)  ,\ldots,\gamma\left(
n\right)  \right)  =\left(  a_{\operatorname*{id}\left(  1\right)
},a_{\operatorname*{id}\left(  2\right)  },\ldots,a_{\operatorname*{id}\left(
k\right)  },b_{\operatorname*{id}\left(  1\right)  },b_{\operatorname*{id}%
\left(  2\right)  },\ldots,b_{\operatorname*{id}\left(  n-k\right)  }\right)
.
\]
Thus,%
\begin{align*}
\left(  \gamma\left(  1\right)  ,\gamma\left(  2\right)  ,\ldots,\gamma\left(
n\right)  \right)   &  =\left(  a_{\operatorname*{id}\left(  1\right)
},a_{\operatorname*{id}\left(  2\right)  },\ldots,a_{\operatorname*{id}\left(
k\right)  },b_{\operatorname*{id}\left(  1\right)  },b_{\operatorname*{id}%
\left(  2\right)  },\ldots,b_{\operatorname*{id}\left(  n-k\right)  }\right)
\\
&  =\left(  a_{1},a_{2},\ldots,a_{k},b_{1},b_{2},\ldots,b_{n-k}\right)  .
\end{align*}
Hence,%
\begin{equation}
\left(  \gamma\left(  i\right)  =a_{i}\ \ \ \ \ \ \ \ \ \ \text{for every
}i\in\left\{  1,2,\ldots,k\right\}  \right)
\label{pf.lem.sol.addexe.jacobi-complement.Ialbe.short.gamma-tuple.pf.gamma1}%
\end{equation}
and%
\begin{equation}
\left(  \gamma\left(  i\right)  =b_{i-k}\ \ \ \ \ \ \ \ \ \ \text{for every
}i\in\left\{  k+1,k+2,\ldots,n\right\}  \right)  .
\label{pf.lem.sol.addexe.jacobi-complement.Ialbe.short.gamma-tuple.pf.gamma2}%
\end{equation}
Now, (\ref{pf.lem.sol.addexe.jacobi-complement.Ialbe.gamma-tuple1}) follows
immediately from
(\ref{pf.lem.sol.addexe.jacobi-complement.Ialbe.short.gamma-tuple.pf.gamma1}).
Furthermore, (\ref{pf.lem.sol.addexe.jacobi-complement.Ialbe.gamma-tuple2})
follows immediately from
(\ref{pf.lem.sol.addexe.jacobi-complement.Ialbe.short.gamma-tuple.pf.gamma2}).
Hence, both (\ref{pf.lem.sol.addexe.jacobi-complement.Ialbe.gamma-tuple1}) and
(\ref{pf.lem.sol.addexe.jacobi-complement.Ialbe.gamma-tuple2}) are proven.]
\end{vershort}

\begin{verlong}
[\textit{Proof of
(\ref{pf.lem.sol.addexe.jacobi-complement.Ialbe.gamma-tuple1}) and
(\ref{pf.lem.sol.addexe.jacobi-complement.Ialbe.gamma-tuple2}):} The
permutation $\sigma_{I,\operatorname*{id}\nolimits_{k},\operatorname*{id}%
\nolimits_{n-k}}$ is the unique $\sigma\in S_{n}$ satisfying%
\begin{equation}
\left(  \sigma\left(  1\right)  ,\sigma\left(  2\right)  ,\ldots,\sigma\left(
n\right)  \right)  =\left(  a_{\operatorname*{id}\nolimits_{k}\left(
1\right)  },a_{\operatorname*{id}\nolimits_{k}\left(  2\right)  }%
,\ldots,a_{\operatorname*{id}\nolimits_{k}\left(  k\right)  }%
,b_{\operatorname*{id}\nolimits_{n-k}\left(  1\right)  },b_{\operatorname*{id}%
\nolimits_{n-k}\left(  2\right)  },\ldots,b_{\operatorname*{id}\nolimits_{n-k}%
\left(  n-k\right)  }\right)
\label{pf.lem.sol.addexe.jacobi-complement.Ialbe.gamma-tuple.pf.1}%
\end{equation}
(because this is how $\sigma_{I,\operatorname*{id}\nolimits_{k}%
,\operatorname*{id}\nolimits_{n-k}}$ is defined). Hence, $\sigma
_{I,\operatorname*{id}\nolimits_{k},\operatorname*{id}\nolimits_{n-k}}$ is a
$\sigma\in S_{n}$ satisfying
(\ref{pf.lem.sol.addexe.jacobi-complement.Ialbe.gamma-tuple.pf.1}). In other
words, $\gamma$ is a $\sigma\in S_{n}$ satisfying
(\ref{pf.lem.sol.addexe.jacobi-complement.Ialbe.gamma-tuple.pf.1}) (since
$\gamma=\sigma_{I,\operatorname*{id}\nolimits_{k},\operatorname*{id}%
\nolimits_{n-k}}$). In other words, $\gamma$ is an element of $S_{n}$ and
satisfies%
\begin{align}
&  \left(  \gamma\left(  1\right)  ,\gamma\left(  2\right)  ,\ldots
,\gamma\left(  n\right)  \right) \nonumber\\
&  =\left(  a_{\operatorname*{id}\nolimits_{k}\left(  1\right)  }%
,a_{\operatorname*{id}\nolimits_{k}\left(  2\right)  },\ldots
,a_{\operatorname*{id}\nolimits_{k}\left(  k\right)  },b_{\operatorname*{id}%
\nolimits_{n-k}\left(  1\right)  },b_{\operatorname*{id}\nolimits_{n-k}\left(
2\right)  },\ldots,b_{\operatorname*{id}\nolimits_{n-k}\left(  n-k\right)
}\right)  . \label{pf.lem.sol.addexe.jacobi-complement.Ialbe.gamma-tuple.pf.2}%
\end{align}
Now, $\left(  \gamma\left(  1\right)  ,\gamma\left(  2\right)  ,\ldots
,\gamma\left(  k\right)  \right)  =\left(  a_{\operatorname*{id}%
\nolimits_{k}\left(  1\right)  },a_{\operatorname*{id}\nolimits_{k}\left(
2\right)  },\ldots,a_{\operatorname*{id}\nolimits_{k}\left(  k\right)
}\right)  $\ \ \ \ \footnote{\textit{Proof.} We have%
\begin{align*}
&  \left(  \gamma\left(  1\right)  ,\gamma\left(  2\right)  ,\ldots
,\gamma\left(  k\right)  \right) \\
&  =\left(  \text{the list of the first }k\text{ entries of the list
}\underbrace{\left(  \gamma\left(  1\right)  ,\gamma\left(  2\right)
,\ldots,\gamma\left(  n\right)  \right)  }_{\substack{=\left(
a_{\operatorname*{id}\nolimits_{k}\left(  1\right)  },a_{\operatorname*{id}%
\nolimits_{k}\left(  2\right)  },\ldots,a_{\operatorname*{id}\nolimits_{k}%
\left(  k\right)  },b_{\operatorname*{id}\nolimits_{n-k}\left(  1\right)
},b_{\operatorname*{id}\nolimits_{n-k}\left(  2\right)  },\ldots
,b_{\operatorname*{id}\nolimits_{n-k}\left(  n-k\right)  }\right)  \\\text{(by
(\ref{pf.lem.sol.addexe.jacobi-complement.Ialbe.gamma-tuple.pf.2}))}}}\right)
\\
&  =\left(  \text{the list of the first }k\text{ entries of the list }\left(
a_{\operatorname*{id}\nolimits_{k}\left(  1\right)  },a_{\operatorname*{id}%
\nolimits_{k}\left(  2\right)  },\ldots,a_{\operatorname*{id}\nolimits_{k}%
\left(  k\right)  },b_{\operatorname*{id}\nolimits_{n-k}\left(  1\right)
},b_{\operatorname*{id}\nolimits_{n-k}\left(  2\right)  },\ldots
,b_{\operatorname*{id}\nolimits_{n-k}\left(  n-k\right)  }\right)  \right) \\
&  =\left(  a_{\operatorname*{id}\nolimits_{k}\left(  1\right)  }%
,a_{\operatorname*{id}\nolimits_{k}\left(  2\right)  },\ldots
,a_{\operatorname*{id}\nolimits_{k}\left(  k\right)  }\right)  ,
\end{align*}
qed.}. In other words, $\gamma\left(  i\right)  =a_{\operatorname*{id}%
\nolimits_{k}\left(  i\right)  }$ for every $i\in\left\{  1,2,\ldots
,k\right\}  $. Hence, for every $i\in\left\{  1,2,\ldots,k\right\}  $, we have
$\gamma\left(  i\right)  =a_{\operatorname*{id}\nolimits_{k}\left(  i\right)
}=a_{i}$ (since $\underbrace{\operatorname*{id}\nolimits_{k}}%
_{=\operatorname*{id}\nolimits_{\left\{  1,2,\ldots,k\right\}  }}\left(
i\right)  =\operatorname*{id}\nolimits_{\left\{  1,2,\ldots,k\right\}
}\left(  i\right)  =i$). In other words, $\left(  \gamma\left(  1\right)
,\gamma\left(  2\right)  ,\ldots,\gamma\left(  k\right)  \right)  =\left(
a_{1},a_{2},\ldots,a_{k}\right)  $. This proves
(\ref{pf.lem.sol.addexe.jacobi-complement.Ialbe.gamma-tuple1}).

Furthermore, $\left(  \gamma\left(  k+1\right)  ,\gamma\left(  k+2\right)
,\ldots,\gamma\left(  n\right)  \right)  =\left(  b_{\operatorname*{id}%
\nolimits_{n-k}\left(  1\right)  },b_{\operatorname*{id}\nolimits_{n-k}\left(
2\right)  },\ldots,b_{\operatorname*{id}\nolimits_{n-k}\left(  n-k\right)
}\right)  $\ \ \ \ \footnote{\textit{Proof.} We have%
\begin{align*}
&  \left(  \gamma\left(  k+1\right)  ,\gamma\left(  k+2\right)  ,\ldots
,\gamma\left(  n\right)  \right) \\
&  =\left(  \text{the list of the last }n-k\text{ entries of the list
}\underbrace{\left(  \gamma\left(  1\right)  ,\gamma\left(  2\right)
,\ldots,\gamma\left(  n\right)  \right)  }_{\substack{=\left(
a_{\operatorname*{id}\nolimits_{k}\left(  1\right)  },a_{\operatorname*{id}%
\nolimits_{k}\left(  2\right)  },\ldots,a_{\operatorname*{id}\nolimits_{k}%
\left(  k\right)  },b_{\operatorname*{id}\nolimits_{n-k}\left(  1\right)
},b_{\operatorname*{id}\nolimits_{n-k}\left(  2\right)  },\ldots
,b_{\operatorname*{id}\nolimits_{n-k}\left(  n-k\right)  }\right)  \\\text{(by
(\ref{pf.lem.sol.addexe.jacobi-complement.Ialbe.gamma-tuple.pf.2}))}}}\right)
\\
&  =\left(  \text{the list of the last }n-k\text{ entries of the list }\left(
a_{\operatorname*{id}\nolimits_{k}\left(  1\right)  },a_{\operatorname*{id}%
\nolimits_{k}\left(  2\right)  },\ldots,a_{\operatorname*{id}\nolimits_{k}%
\left(  k\right)  },b_{\operatorname*{id}\nolimits_{n-k}\left(  1\right)
},b_{\operatorname*{id}\nolimits_{n-k}\left(  2\right)  },\ldots
,b_{\operatorname*{id}\nolimits_{n-k}\left(  n-k\right)  }\right)  \right) \\
&  =\left(  b_{\operatorname*{id}\nolimits_{n-k}\left(  1\right)
},b_{\operatorname*{id}\nolimits_{n-k}\left(  2\right)  },\ldots
,b_{\operatorname*{id}\nolimits_{n-k}\left(  n-k\right)  }\right)  ,
\end{align*}
qed.}. In other words, $\gamma\left(  k+i\right)  =b_{\operatorname*{id}%
\nolimits_{n-k}\left(  i\right)  }$ for every $i\in\left\{  1,2,\ldots
,n-k\right\}  $. Hence, for every $i\in\left\{  1,2,\ldots,n-k\right\}  $, we
have $\gamma\left(  k+i\right)  =b_{\operatorname*{id}\nolimits_{n-k}\left(
i\right)  }=b_{i}$ (since $\underbrace{\operatorname*{id}\nolimits_{n-k}%
}_{=\operatorname*{id}\nolimits_{\left\{  1,2,\ldots,n-k\right\}  }}\left(
i\right)  =\operatorname*{id}\nolimits_{\left\{  1,2,\ldots,n-k\right\}
}\left(  i\right)  =i$). In other words, $\left(  \gamma\left(  k+1\right)
,\gamma\left(  k+2\right)  ,\ldots,\gamma\left(  n\right)  \right)  =\left(
b_{1},b_{2},\ldots,b_{n-k}\right)  $. This proves
(\ref{pf.lem.sol.addexe.jacobi-complement.Ialbe.gamma-tuple2}). Hence, both
(\ref{pf.lem.sol.addexe.jacobi-complement.Ialbe.gamma-tuple1}) and
(\ref{pf.lem.sol.addexe.jacobi-complement.Ialbe.gamma-tuple2}) are proven.]
\end{verlong}

Comparing (\ref{pf.lem.sol.addexe.jacobi-complement.Ialbe.gamma-tuple1}) with
$w\left(  I\right)  =\left(  a_{1},a_{2},\ldots,a_{k}\right)  $, we obtain
\[
\left(  \gamma\left(  1\right)  ,\gamma\left(  2\right)  ,\ldots,\gamma\left(
k\right)  \right)  =w\left(  I\right)  .
\]
Comparing (\ref{pf.lem.sol.addexe.jacobi-complement.Ialbe.gamma-tuple2}) with
$w\left(  \widetilde{I}\right)  =\left(  b_{1},b_{2},\ldots,b_{n-k}\right)  $,
we obtain
\[
\left(  \gamma\left(  k+1\right)  ,\gamma\left(  k+2\right)  ,\ldots
,\gamma\left(  n\right)  \right)  =w\left(  \widetilde{I}\right)  .
\]
Now, we have shown that the permutation $\gamma\in S_{n}$ satisfies $\left(
\gamma\left(  1\right)  ,\gamma\left(  2\right)  ,\ldots,\gamma\left(
k\right)  \right)  =w\left(  I\right)  $, $\left(  \gamma\left(  k+1\right)
,\gamma\left(  k+2\right)  ,\ldots,\gamma\left(  n\right)  \right)  =w\left(
\widetilde{I}\right)  $ and $\left(  -1\right)  ^{\gamma}=\left(  -1\right)
^{\sum I-\left(  1+2+\cdots+k\right)  }$. Hence, there exists a $\sigma\in
S_{n}$ satisfying $\left(  \sigma\left(  1\right)  ,\sigma\left(  2\right)
,\ldots,\sigma\left(  k\right)  \right)  =w\left(  I\right)  $, \newline%
$\left(  \sigma\left(  k+1\right)  ,\sigma\left(  k+2\right)  ,\ldots
,\sigma\left(  n\right)  \right)  =w\left(  \widetilde{I}\right)  $ and
$\left(  -1\right)  ^{\sigma}=\left(  -1\right)  ^{\sum I-\left(
1+2+\cdots+k\right)  }$ (namely, $\sigma=\gamma$). This proves Lemma
\ref{lem.sol.addexe.jacobi-complement.Ialbe}.
\end{proof}

\begin{proof}
[Proof of Lemma \ref{lem.det.laplace-multi.Apq}.]Let us first study the sets
$P$ and $\widetilde{P}$.

\begin{vershort}
Define $k\in\mathbb{N}$ by $k=\left\vert P\right\vert =\left\vert Q\right\vert
$. (This makes sense, since we assumed that $\left\vert P\right\vert
=\left\vert Q\right\vert $.)

The definition of $\widetilde{P}$ yields $\widetilde{P}=\left\{
1,2,\ldots,n\right\}  \setminus P$. Since $P$ is a subset of $\left\{
1,2,\ldots,n\right\}  $, this leads to $\left\vert \widetilde{P}\right\vert
=\underbrace{\left\vert \left\{  1,2,\ldots,n\right\}  \right\vert }%
_{=n}-\underbrace{\left\vert P\right\vert }_{=k}=n-k$. Thus, $n-k=\left\vert
\widetilde{P}\right\vert \geq0$, so that $n\geq k$ and thus $k\in\left\{
0,1,\ldots,n\right\}  $.
\end{vershort}

\begin{verlong}
The definition of $\widetilde{P}$ yields $\widetilde{P}=\left\{
1,2,\ldots,n\right\}  \setminus P\subseteq\left\{  1,2,\ldots,n\right\}  $.
Thus, $\widetilde{P}$ is a finite set (since $\left\{  1,2,\ldots,n\right\}  $
is a finite set). Hence, $\left\vert \widetilde{P}\right\vert \in\mathbb{N}$.

Also, $P\subseteq\left\{  1,2,\ldots,n\right\}  $. Thus, $P$ is a finite set
(since $\left\{  1,2,\ldots,n\right\}  $ is a finite set). Hence, $\left\vert
P\right\vert \in\mathbb{N}$. Thus, we can define $k\in\mathbb{N}$ by
$k=\left\vert P\right\vert $. Consider this $k$. We have $k=\left\vert
P\right\vert =\left\vert Q\right\vert $.

Also,%
\begin{align*}
\left\vert \underbrace{\widetilde{P}}_{=\left\{  1,2,\ldots,n\right\}
\setminus P}\right\vert  &  =\left\vert \left\{  1,2,\ldots,n\right\}
\setminus P\right\vert =\underbrace{\left\vert \left\{  1,2,\ldots,n\right\}
\right\vert }_{=n}-\underbrace{\left\vert P\right\vert }_{=k}%
\ \ \ \ \ \ \ \ \ \ \left(  \text{since }P\subseteq\left\{  1,2,\ldots
,n\right\}  \right) \\
&  =n-k.
\end{align*}
Hence, $n-k=\left\vert \widetilde{P}\right\vert \in\mathbb{N}$; thus,
$n-k\geq0$. In other words, $n\geq k$. In other words, $k\leq n$. Combined
with $k\geq0$ (since $k=\left\vert P\right\vert \in\mathbb{N}$), this yields
$k\in\left\{  0,1,\ldots,n\right\}  $.
\end{verlong}

Let $\left(  p_{1},p_{2},\ldots,p_{k}\right)  $ be the list of all elements of
$P$ in increasing order (with no repetitions). (This is well-defined, because
$\left\vert P\right\vert =k$.)

Let $\left(  r_{1},r_{2},\ldots,r_{n-k}\right)  $ be the list of all elements
of $\widetilde{P}$ in increasing order (with no repetitions). (This is
well-defined, because $\left\vert \widetilde{P}\right\vert =n-k$.)

Lemma \ref{lem.sol.addexe.jacobi-complement.Ialbe} (applied to $I=P$) shows
that there exists a $\sigma\in S_{n}$ satisfying $\left(  \sigma\left(
1\right)  ,\sigma\left(  2\right)  ,\ldots,\sigma\left(  k\right)  \right)
=w\left(  P\right)  $, $\left(  \sigma\left(  k+1\right)  ,\sigma\left(
k+2\right)  ,\ldots,\sigma\left(  n\right)  \right)  =w\left(  \widetilde{P}%
\right)  $ and $\left(  -1\right)  ^{\sigma}=\left(  -1\right)  ^{\sum
P-\left(  1+2+\cdots+k\right)  }$. Denote this $\sigma$ by $\gamma$. Thus,
$\gamma$ is an element of $S_{n}$ satisfying $\left(  \gamma\left(  1\right)
,\gamma\left(  2\right)  ,\ldots,\gamma\left(  k\right)  \right)  =w\left(
P\right)  $, $\left(  \gamma\left(  k+1\right)  ,\gamma\left(  k+2\right)
,\ldots,\gamma\left(  n\right)  \right)  =w\left(  \widetilde{P}\right)  $ and
$\left(  -1\right)  ^{\gamma}=\left(  -1\right)  ^{\sum P-\left(
1+2+\cdots+k\right)  }$.

Notice that%
\begin{align}
\left(  -1\right)  ^{\sum P-\left(  1+2+\cdots+k\right)  }\underbrace{\left(
-1\right)  ^{\gamma}}_{=\left(  -1\right)  ^{\sum P-\left(  1+2+\cdots
+k\right)  }}  &  =\left(  -1\right)  ^{\sum P-\left(  1+2+\cdots+k\right)
}\left(  -1\right)  ^{\sum P-\left(  1+2+\cdots+k\right)  }\nonumber\\
&  =\left(  \left(  -1\right)  ^{\sum P-\left(  1+2+\cdots+k\right)  }\right)
^{2}=\left(  -1\right)  ^{2\left(  \sum P-\left(  1+2+\cdots+k\right)
\right)  }\nonumber\\
&  =1 \label{pf.lem.det.laplace-multi.Apq.sign-gamma-multed}%
\end{align}
(since $2\left(  \sum P-\left(  1+2+\cdots+k\right)  \right)  $ is even).

Now, we claim that%
\begin{equation}
\left(  \gamma\left(  i\right)  =p_{i}\ \ \ \ \ \ \ \ \ \ \text{for every
}i\in\left\{  1,2,\ldots,k\right\}  \right)
\label{pf.lem.det.laplace-multi.Apq.gamma-tuple1}%
\end{equation}
and%
\begin{equation}
\left(  \gamma\left(  k+i\right)  =r_{i}\ \ \ \ \ \ \ \ \ \ \text{for every
}i\in\left\{  1,2,\ldots,n-k\right\}  \right)  .
\label{pf.lem.det.laplace-multi.Apq.gamma-tuple2}%
\end{equation}

\begin{vershort}
[\textit{Proof of (\ref{pf.lem.det.laplace-multi.Apq.gamma-tuple1}) and
(\ref{pf.lem.det.laplace-multi.Apq.gamma-tuple2}):} The lists $w\left(
P\right)  $ and $\left(  p_{1},p_{2},\ldots,p_{k}\right)  $ must be identical
(since they are both defined to be the list of all elements of $P$ in
increasing order (with no repetitions)). In other words, we have $w\left(
P\right)  =\left(  p_{1},p_{2},\ldots,p_{k}\right)  $. Now,
\[
\left(  \gamma\left(  1\right)  ,\gamma\left(  2\right)  ,\ldots,\gamma\left(
k\right)  \right)  =w\left(  P\right)  =\left(  p_{1},p_{2},\ldots
,p_{k}\right)  .
\]
In other words, $\gamma\left(  i\right)  =p_{i}$ for every $i\in\left\{
1,2,\ldots,k\right\}  $. This proves
(\ref{pf.lem.det.laplace-multi.Apq.gamma-tuple1}).

The lists $w\left(  \widetilde{P}\right)  $ and $\left(  r_{1},r_{2}%
,\ldots,r_{n-k}\right)  $ must be identical (since they are both defined to be
the list of all elements of $\widetilde{P}$ in increasing order (with no
repetitions)). In other words, we have $w\left(  \widetilde{P}\right)
=\left(  r_{1},r_{2},\ldots,r_{n-k}\right)  $. Now,%
\[
\left(  \gamma\left(  k+1\right)  ,\gamma\left(  k+2\right)  ,\ldots
,\gamma\left(  n\right)  \right)  =w\left(  \widetilde{P}\right)  =\left(
r_{1},r_{2},\ldots,r_{n-k}\right)  .
\]
In other words, $\gamma\left(  k+i\right)  =r_{i}$ for every $i\in\left\{
1,2,\ldots,n-k\right\}  $. This proves
(\ref{pf.lem.det.laplace-multi.Apq.gamma-tuple2}). Thus, both
(\ref{pf.lem.det.laplace-multi.Apq.gamma-tuple1}) and
(\ref{pf.lem.det.laplace-multi.Apq.gamma-tuple2}) are proven.]
\end{vershort}

\begin{verlong}
[\textit{Proof of (\ref{pf.lem.det.laplace-multi.Apq.gamma-tuple1}) and
(\ref{pf.lem.det.laplace-multi.Apq.gamma-tuple2}):} Recall that $w\left(
P\right)  $ is the list of all elements of $P$ in increasing order (with no
repetitions) (by the definition of $w\left(  P\right)  $). Thus,%
\begin{align*}
w\left(  P\right)   &  =\left(  \text{the list of all elements of }P\text{ in
increasing order (with no repetitions)}\right) \\
&  =\left(  p_{1},p_{2},\ldots,p_{k}\right)
\end{align*}
(since $\left(  p_{1},p_{2},\ldots,p_{k}\right)  $ is the list of all elements
of $P$ in increasing order (with no repetitions)). Now,
\[
\left(  \gamma\left(  1\right)  ,\gamma\left(  2\right)  ,\ldots,\gamma\left(
k\right)  \right)  =w\left(  P\right)  =\left(  p_{1},p_{2},\ldots
,p_{k}\right)  .
\]
In other words, $\gamma\left(  i\right)  =p_{i}$ for every $i\in\left\{
1,2,\ldots,k\right\}  $. This proves
(\ref{pf.lem.det.laplace-multi.Apq.gamma-tuple1}).

Recall that $w\left(  \widetilde{P}\right)  $ is the list of all elements of
$\widetilde{P}$ in increasing order (with no repetitions) (by the definition
of $w\left(  \widetilde{P}\right)  $). Thus,%
\begin{align*}
w\left(  \widetilde{P}\right)   &  =\left(  \text{the list of all elements of
}\widetilde{P}\text{ in increasing order (with no repetitions)}\right) \\
&  =\left(  r_{1},r_{2},\ldots,r_{n-k}\right)
\end{align*}
(since $\left(  r_{1},r_{2},\ldots,r_{n-k}\right)  $ is the list of all
elements of $\widetilde{P}$ in increasing order (with no repetitions)). Now,%
\[
\left(  \gamma\left(  k+1\right)  ,\gamma\left(  k+2\right)  ,\ldots
,\gamma\left(  n\right)  \right)  =w\left(  \widetilde{P}\right)  =\left(
r_{1},r_{2},\ldots,r_{n-k}\right)  .
\]
In other words, $\gamma\left(  k+i\right)  =r_{i}$ for every $i\in\left\{
1,2,\ldots,n-k\right\}  $. This proves
(\ref{pf.lem.det.laplace-multi.Apq.gamma-tuple2}). Thus, both
(\ref{pf.lem.det.laplace-multi.Apq.gamma-tuple1}) and
(\ref{pf.lem.det.laplace-multi.Apq.gamma-tuple2}) are proven.]
\end{verlong}

\begin{verlong}
Now, recall that $\gamma$ is an element of $S_{n}$. In other words, $\gamma$
is a permutation of $\left\{  1,2,\ldots,n\right\}  $ (since $S_{n}$ is the
set of all permutations of $\left\{  1,2,\ldots,n\right\}  $). Hence, $\gamma$
is a bijective map $\left\{  1,2,\ldots,n\right\}  \rightarrow\left\{
1,2,\ldots,n\right\}  $.
\end{verlong}

Define an $n\times n$-matrix $A^{\prime}$ by $A^{\prime}=\left(
a_{\gamma\left(  i\right)  ,j}\right)  _{1\leq i\leq n,\ 1\leq j\leq n}$.
Define an $n\times n$-matrix $B^{\prime}$ by $B^{\prime}=\left(
b_{\gamma\left(  i\right)  ,j}\right)  _{1\leq i\leq n,\ 1\leq j\leq n}$.
Then,%
\begin{equation}
\operatorname*{sub}\nolimits_{w\left(  P\right)  }^{w\left(  Q\right)
}A=\operatorname*{sub}\nolimits_{\left(  1,2,\ldots,k\right)  }^{w\left(
Q\right)  }\left(  A^{\prime}\right)
\label{pf.lem.det.laplace-multi.Apq.sub1}%
\end{equation}
\footnote{\textit{Proof of (\ref{pf.lem.det.laplace-multi.Apq.sub1}):} Write
the list $w\left(  Q\right)  $ in the form $w\left(  Q\right)  =\left(
q_{1},q_{2},\ldots,q_{\ell}\right)  $ for some $\ell\in\mathbb{N}$. (Actually,
$\ell=k$, but we will not use this.)
\par
Recall that $\left(  \gamma\left(  1\right)  ,\gamma\left(  2\right)
,\ldots,\gamma\left(  k\right)  \right)  =w\left(  P\right)  $, so that
$w\left(  P\right)  =\left(  \gamma\left(  1\right)  ,\gamma\left(  2\right)
,\ldots,\gamma\left(  k\right)  \right)  $.
\par
From $w\left(  P\right)  =\left(  \gamma\left(  1\right)  ,\gamma\left(
2\right)  ,\ldots,\gamma\left(  k\right)  \right)  $ and $w\left(  Q\right)
=\left(  q_{1},q_{2},\ldots,q_{\ell}\right)  $, we obtain%
\begin{align}
\operatorname*{sub}\nolimits_{w\left(  P\right)  }^{w\left(  Q\right)  }A  &
=\operatorname*{sub}\nolimits_{\left(  \gamma\left(  1\right)  ,\gamma\left(
2\right)  ,\ldots,\gamma\left(  k\right)  \right)  }^{\left(  q_{1}%
,q_{2},\ldots,q_{\ell}\right)  }A=\operatorname*{sub}\nolimits_{\gamma\left(
1\right)  ,\gamma\left(  2\right)  ,\ldots,\gamma\left(  k\right)  }%
^{q_{1},q_{2},\ldots,q_{\ell}}A=\left(  a_{\gamma\left(  x\right)  ,q_{y}%
}\right)  _{1\leq x\leq k,\ 1\leq y\leq\ell}%
\label{pf.lem.det.laplace-multi.Apq.sub1.pf.1}\\
&  \ \ \ \ \ \ \ \ \ \ \left(  \text{by the definition of }\operatorname*{sub}%
\nolimits_{\gamma\left(  1\right)  ,\gamma\left(  2\right)  ,\ldots
,\gamma\left(  k\right)  }^{q_{1},q_{2},\ldots,q_{\ell}}A\text{, since
}A=\left(  a_{i,j}\right)  _{1\leq i\leq n,\ 1\leq j\leq n}\right)  .\nonumber
\end{align}
\par
On the other hand, from $w\left(  Q\right)  =\left(  q_{1},q_{2}%
,\ldots,q_{\ell}\right)  $, we obtain%
\[
\operatorname*{sub}\nolimits_{\left(  1,2,\ldots,k\right)  }^{w\left(
Q\right)  }\left(  A^{\prime}\right)  =\operatorname*{sub}\nolimits_{\left(
1,2,\ldots,k\right)  }^{\left(  q_{1},q_{2},\ldots,q_{\ell}\right)  }\left(
A^{\prime}\right)  =\operatorname*{sub}\nolimits_{1,2,\ldots,k}^{q_{1}%
,q_{2},\ldots,q_{\ell}}\left(  A^{\prime}\right)  =\left(  a_{\gamma\left(
x\right)  ,q_{y}}\right)  _{1\leq x\leq k,\ 1\leq y\leq\ell}%
\]
(by the definition of $\operatorname*{sub}\nolimits_{1,2,\ldots,k}%
^{q_{1},q_{2},\ldots,q_{\ell}}\left(  A^{\prime}\right)  $, since $A^{\prime
}=\left(  a_{\gamma\left(  i\right)  ,j}\right)  _{1\leq i\leq n,\ 1\leq j\leq
n}$). Comparing this with (\ref{pf.lem.det.laplace-multi.Apq.sub1.pf.1}), we
obtain $\operatorname*{sub}\nolimits_{w\left(  P\right)  }^{w\left(  Q\right)
}A=\operatorname*{sub}\nolimits_{\left(  1,2,\ldots,k\right)  }^{w\left(
Q\right)  }\left(  A^{\prime}\right)  $. This proves
(\ref{pf.lem.det.laplace-multi.Apq.sub1}).} and%
\begin{equation}
\operatorname*{sub}\nolimits_{w\left(  \widetilde{P}\right)  }^{w\left(
\widetilde{Q}\right)  }B=\operatorname*{sub}\nolimits_{\left(  k+1,k+2,\ldots
,n\right)  }^{w\left(  \widetilde{Q}\right)  }\left(  B^{\prime}\right)
\label{pf.lem.det.laplace-multi.Apq.sub2}%
\end{equation}
\footnote{\textit{Proof of (\ref{pf.lem.det.laplace-multi.Apq.sub2}):} Write
the list $w\left(  \widetilde{Q}\right)  $ in the form $w\left(
\widetilde{Q}\right)  =\left(  q_{1},q_{2},\ldots,q_{\ell}\right)  $ for some
$\ell\in\mathbb{N}$. (Actually, $\ell=n-k$, but we will not use this.)
\par
Recall that $\left(  \gamma\left(  k+1\right)  ,\gamma\left(  k+2\right)
,\ldots,\gamma\left(  n\right)  \right)  =w\left(  \widetilde{P}\right)  $.
Thus, $w\left(  \widetilde{P}\right)  =\left(  \gamma\left(  k+1\right)
,\gamma\left(  k+2\right)  ,\ldots,\gamma\left(  n\right)  \right)  $.
\par
From $w\left(  \widetilde{P}\right)  =\left(  \gamma\left(  k+1\right)
,\gamma\left(  k+2\right)  ,\ldots,\gamma\left(  n\right)  \right)  $ and
$w\left(  \widetilde{Q}\right)  =\left(  q_{1},q_{2},\ldots,q_{\ell}\right)
$, we obtain%
\begin{align}
\operatorname*{sub}\nolimits_{w\left(  \widetilde{P}\right)  }^{w\left(
\widetilde{Q}\right)  }B  &  =\operatorname*{sub}\nolimits_{\left(
\gamma\left(  k+1\right)  ,\gamma\left(  k+2\right)  ,\ldots,\gamma\left(
n\right)  \right)  }^{\left(  q_{1},q_{2},\ldots,q_{\ell}\right)
}B=\operatorname*{sub}\nolimits_{\gamma\left(  k+1\right)  ,\gamma\left(
k+2\right)  ,\ldots,\gamma\left(  n\right)  }^{q_{1},q_{2},\ldots,q_{\ell}%
}B=\left(  b_{\gamma\left(  k+x\right)  ,q_{y}}\right)  _{1\leq x\leq
n-k,\ 1\leq y\leq\ell}\label{pf.lem.det.laplace-multi.Apq.sub2.pf.1}\\
&  \ \ \ \ \ \ \ \ \ \ \left(  \text{by the definition of }\operatorname*{sub}%
\nolimits_{\gamma\left(  k+1\right)  ,\gamma\left(  k+2\right)  ,\ldots
,\gamma\left(  n\right)  }^{q_{1},q_{2},\ldots,q_{\ell}}B\text{, since
}B=\left(  b_{i,j}\right)  _{1\leq i\leq n,\ 1\leq j\leq n}\right)  .\nonumber
\end{align}
\par
On the other hand, from $w\left(  \widetilde{Q}\right)  =\left(  q_{1}%
,q_{2},\ldots,q_{\ell}\right)  $, we obtain%
\[
\operatorname*{sub}\nolimits_{\left(  k+1,k+2,\ldots,n\right)  }^{w\left(
\widetilde{Q}\right)  }\left(  B^{\prime}\right)  =\operatorname*{sub}%
\nolimits_{\left(  k+1,k+2,\ldots,n\right)  }^{\left(  q_{1},q_{2}%
,\ldots,q_{\ell}\right)  }\left(  B^{\prime}\right)  =\operatorname*{sub}%
\nolimits_{k+1,k+2,\ldots,n}^{q_{1},q_{2},\ldots,q_{\ell}}\left(  B^{\prime
}\right)  =\left(  b_{\gamma\left(  k+x\right)  ,q_{y}}\right)  _{1\leq x\leq
n-k,\ 1\leq y\leq\ell}%
\]
(by the definition of $\operatorname*{sub}\nolimits_{k+1,k+2,\ldots,n}%
^{q_{1},q_{2},\ldots,q_{\ell}}\left(  B^{\prime}\right)  $, since $B^{\prime
}=\left(  b_{\gamma\left(  i\right)  ,j}\right)  _{1\leq i\leq n,\ 1\leq j\leq
n}$). Comparing this with (\ref{pf.lem.det.laplace-multi.Apq.sub2.pf.1}), we
obtain $\operatorname*{sub}\nolimits_{w\left(  \widetilde{P}\right)
}^{w\left(  \widetilde{Q}\right)  }B=\operatorname*{sub}\nolimits_{\left(
k+1,k+2,\ldots,n\right)  }^{w\left(  \widetilde{Q}\right)  }\left(  B^{\prime
}\right)  $. This proves (\ref{pf.lem.det.laplace-multi.Apq.sub2}).}.

Lemma \ref{lem.sol.det.laplace-multi.group-bij} shows that the map
$S_{n}\rightarrow S_{n},\ \sigma\mapsto\sigma\circ\gamma$ is a bijection.

But $\gamma\left(  \left\{  1,2,\ldots,k\right\}  \right)  =P$%
\ \ \ \ \footnote{\textit{Proof.} We have
\begin{align*}
\gamma\left(  \left\{  1,2,\ldots,k\right\}  \right)   &  =\left\{
\gamma\left(  1\right)  ,\gamma\left(  2\right)  ,\ldots,\gamma\left(
k\right)  \right\}  =\left\{  \underbrace{\gamma\left(  i\right)
}_{\substack{=p_{i}\\\text{(by
(\ref{pf.lem.det.laplace-multi.Apq.gamma-tuple1}))}}}\ \mid\ i\in\left\{
1,2,\ldots,k\right\}  \right\} \\
&  =\left\{  p_{i}\ \mid\ i\in\left\{  1,2,\ldots,k\right\}  \right\}
=\left\{  p_{1},p_{2},\ldots,p_{k}\right\}  =P
\end{align*}
(since $\left(  p_{1},p_{2},\ldots,p_{k}\right)  $ is a list of all elements
of $P$). Qed.}. Hence, every $\sigma\in S_{n}$ satisfies%
\begin{equation}
\left(  \sigma\circ\gamma\right)  \left(  \left\{  1,2,\ldots,k\right\}
\right)  =\sigma\left(  \underbrace{\gamma\left(  \left\{  1,2,\ldots
,k\right\}  \right)  }_{=P}\right)  =\sigma\left(  P\right)  .
\label{pf.lem.det.laplace-multi.Apq.3a}%
\end{equation}

\begin{vershort}
Furthermore, every $\sigma\in S_{n}$ satisfies%
\begin{equation}
\prod_{i\in\left\{  1,2,\ldots,k\right\}  }a_{\gamma\left(  i\right)  ,\left(
\sigma\circ\gamma\right)  \left(  i\right)  }=\prod_{i\in P}a_{i,\sigma\left(
i\right)  } \label{pf.lem.det.laplace-multi.Apq.short.3b}%
\end{equation}
\footnote{\textit{Proof of (\ref{pf.lem.det.laplace-multi.Apq.short.3b}):} Let
$\sigma\in S_{n}$.
\par
We shall use the notation introduced in Definition \ref{def.sol.Ialbe.12n}.
Thus, $\left[  k\right]  =\left\{  1,2,\ldots,k\right\}  $.
\par
Every $i\in\left\{  1,2,\ldots,k\right\}  $ satisfies $p_{i}=\gamma\left(
i\right)  $ (by (\ref{pf.lem.det.laplace-multi.Apq.gamma-tuple1})) and
$\sigma\left(  \underbrace{p_{i}}_{=\gamma\left(  i\right)  }\right)
=\sigma\left(  \gamma\left(  i\right)  \right)  =\left(  \sigma\circ
\gamma\right)  \left(  i\right)  $. Hence, every $i\in\left\{  1,2,\ldots
,k\right\}  $ satisfies%
\begin{equation}
a_{p_{i},\sigma\left(  p_{i}\right)  }=a_{\gamma\left(  i\right)  ,\left(
\sigma\circ\gamma\right)  \left(  i\right)  }.
\label{pf.lem.det.laplace-multi.Apq.short.3b.pf.1c}%
\end{equation}
\par
Recall that $\left(  p_{1},p_{2},\ldots,p_{k}\right)  $ is the list of all
elements of $P$ in increasing order (with no repetitions). Hence, Lemma
\ref{lem.sol.Ialbe.inclist} \textbf{(a)} (applied to $P$, $k$ and $\left(
p_{1},p_{2},\ldots,p_{k}\right)  $ instead of $S$, $s$ and $\left(
c_{1},c_{2},\ldots,c_{s}\right)  $) shows that the map $\left[  k\right]
\rightarrow P,\ h\mapsto p_{h}$ is well-defined and a bijection. Hence, we can
substitute $p_{h}$ for $i$ in the product $\prod_{i\in P}a_{i,\sigma\left(
i\right)  }$. We thus obtain%
\begin{align*}
\prod_{i\in P}a_{i,\sigma\left(  i\right)  }  &  =\prod_{h\in\left[  k\right]
}a_{p_{h},\sigma\left(  p_{h}\right)  }=\underbrace{\prod_{i\in\left[
k\right]  }}_{\substack{=\prod_{i\in\left\{  1,2,\ldots,k\right\}
}\\\text{(since }\left[  k\right]  =\left\{  1,2,\ldots,k\right\}  \text{)}%
}}a_{p_{i},\sigma\left(  p_{i}\right)  }\\
&  \ \ \ \ \ \ \ \ \ \ \left(  \text{here, we have renamed the index }h\text{
as }i\text{ in the product}\right) \\
&  =\prod_{i\in\left\{  1,2,\ldots,k\right\}  }\underbrace{a_{p_{i}%
,\sigma\left(  p_{i}\right)  }}_{\substack{=a_{\gamma\left(  i\right)
,\left(  \sigma\circ\gamma\right)  \left(  i\right)  }\\\text{(by
(\ref{pf.lem.det.laplace-multi.Apq.short.3b.pf.1c}))}}}=\prod_{i\in\left\{
1,2,\ldots,k\right\}  }a_{\gamma\left(  i\right)  ,\left(  \sigma\circ
\gamma\right)  \left(  i\right)  }.
\end{align*}
This proves (\ref{pf.lem.det.laplace-multi.Apq.short.3b}).} and%
\begin{equation}
\prod_{i\in\left\{  k+1,k+2,\ldots,n\right\}  }b_{\gamma\left(  i\right)
,\left(  \sigma\circ\gamma\right)  \left(  i\right)  }=\prod_{i\in
\widetilde{P}}b_{i,\sigma\left(  i\right)  }
\label{pf.lem.det.laplace-multi.Apq.short.3c}%
\end{equation}
\footnote{\textit{Proof of (\ref{pf.lem.det.laplace-multi.Apq.short.3c}):} Let
$\sigma\in S_{n}$.
\par
We shall use the notation introduced in Definition \ref{def.sol.Ialbe.12n}.
Thus, $\left[  n-k\right]  =\left\{  1,2,\ldots,n-k\right\}  $.
\par
Every $i\in\left\{  1,2,\ldots,n-k\right\}  $ satisfies $r_{i}=\gamma\left(
k+i\right)  $ (by (\ref{pf.lem.det.laplace-multi.Apq.gamma-tuple2})) and
$\sigma\left(  \underbrace{r_{i}}_{=\gamma\left(  k+i\right)  }\right)
=\sigma\left(  \gamma\left(  k+i\right)  \right)  =\left(  \sigma\circ
\gamma\right)  \left(  k+i\right)  $. Hence, every $i\in\left\{
1,2,\ldots,n-k\right\}  $ satisfies%
\begin{equation}
b_{r_{i},\sigma\left(  r_{i}\right)  }=b_{\gamma\left(  k+i\right)  ,\left(
\sigma\circ\gamma\right)  \left(  k+i\right)  }.
\label{pf.lem.det.laplace-multi.Apq.short.3c.pf.1c}%
\end{equation}
\par
Recall that $\left(  r_{1},r_{2},\ldots,r_{n-k}\right)  $ is the list of all
elements of $\widetilde{P}$ in increasing order (with no repetitions). Hence,
Lemma \ref{lem.sol.Ialbe.inclist} \textbf{(a)} (applied to $\widetilde{P}$,
$n-k$ and $\left(  r_{1},r_{2},\ldots,r_{n-k}\right)  $ instead of $S$, $s$
and $\left(  c_{1},c_{2},\ldots,c_{s}\right)  $) shows that the map $\left[
n-k\right]  \rightarrow\widetilde{P},\ h\mapsto r_{h}$ is well-defined and a
bijection. Hence, we can substitute $r_{h}$ for $i$ in the product
$\prod_{i\in\widetilde{P}}b_{i,\sigma\left(  i\right)  }$. We thus obtain%
\begin{align*}
\prod_{i\in\widetilde{P}}b_{i,\sigma\left(  i\right)  }  &  =\prod
_{h\in\left[  n-k\right]  }b_{r_{h},\sigma\left(  r_{h}\right)  }%
=\underbrace{\prod_{i\in\left[  n-k\right]  }}_{\substack{=\prod_{i\in\left\{
1,2,\ldots,n-k\right\}  }\\\text{(since }\left[  n-k\right]  =\left\{
1,2,\ldots,n-k\right\}  \text{)}}}b_{r_{i},\sigma\left(  r_{i}\right)  }\\
&  \ \ \ \ \ \ \ \ \ \ \left(  \text{here, we have renamed the index }h\text{
as }i\text{ in the product}\right) \\
&  =\prod_{i\in\left\{  1,2,\ldots,n-k\right\}  }\underbrace{b_{r_{i}%
,\sigma\left(  r_{i}\right)  }}_{\substack{=b_{\gamma\left(  k+i\right)
,\left(  \sigma\circ\gamma\right)  \left(  k+i\right)  }\\\text{(by
(\ref{pf.lem.det.laplace-multi.Apq.short.3c.pf.1c}))}}}\\
&  =\prod_{i\in\left\{  1,2,\ldots,n-k\right\}  }b_{\gamma\left(  k+i\right)
,\left(  \sigma\circ\gamma\right)  \left(  k+i\right)  }=\prod_{i\in\left\{
k+1,k+2,\ldots,n\right\}  }b_{\gamma\left(  i\right)  ,\left(  \sigma
\circ\gamma\right)  \left(  i\right)  }.\\
&  \ \ \ \ \ \ \ \ \ \ \left(  \text{here, we have substituted }i\text{ for
}k+i\text{ in the product}\right)  .
\end{align*}
This proves (\ref{pf.lem.det.laplace-multi.Apq.short.3c}).}.
\end{vershort}

\begin{verlong}
Furthermore, every $\sigma\in S_{n}$ satisfies%
\begin{equation}
\prod_{i\in\left\{  1,2,\ldots,k\right\}  }a_{\gamma\left(  i\right)  ,\left(
\sigma\circ\gamma\right)  \left(  i\right)  }=\prod_{i\in P}a_{i,\sigma\left(
i\right)  } \label{pf.lem.det.laplace-multi.Apq.3b}%
\end{equation}
\footnote{\textit{Proof of (\ref{pf.lem.det.laplace-multi.Apq.3b}):} Let
$\sigma\in S_{n}$.
\par
We shall use the notation introduced in Definition \ref{def.sol.Ialbe.12n}.
Thus, $\left[  k\right]  =\left\{  1,2,\ldots,k\right\}  $ (by the definition
of $\left[  k\right]  $).
\par
Every $i\in\left\{  1,2,\ldots,k\right\}  $ satisfies%
\begin{equation}
p_{i}=\gamma\left(  i\right)  \ \ \ \ \ \ \ \ \ \ \left(  \text{by
(\ref{pf.lem.det.laplace-multi.Apq.gamma-tuple1})}\right)
\label{pf.lem.det.laplace-multi.Apq.3b.pf.1a}%
\end{equation}
and%
\begin{equation}
\sigma\left(  \underbrace{p_{i}}_{=\gamma\left(  i\right)  }\right)
=\sigma\left(  \gamma\left(  i\right)  \right)  =\left(  \sigma\circ
\gamma\right)  \left(  i\right)  .
\label{pf.lem.det.laplace-multi.Apq.3b.pf.1b}%
\end{equation}
Hence, every $i\in\left\{  1,2,\ldots,k\right\}  $ satisfies%
\begin{align}
a_{p_{i},\sigma\left(  p_{i}\right)  }  &  =a_{p_{i},\left(  \sigma\circ
\gamma\right)  \left(  i\right)  }\ \ \ \ \ \ \ \ \ \ \left(  \text{since
}\sigma\left(  p_{i}\right)  =\left(  \sigma\circ\gamma\right)  \left(
i\right)  \text{ (by (\ref{pf.lem.det.laplace-multi.Apq.3b.pf.1b}))}\right)
\nonumber\\
&  =a_{\gamma\left(  i\right)  ,\left(  \sigma\circ\gamma\right)  \left(
i\right)  }\ \ \ \ \ \ \ \ \ \ \left(  \text{since }p_{i}=\gamma\left(
i\right)  \text{ (by (\ref{pf.lem.det.laplace-multi.Apq.3b.pf.1a}))}\right)  .
\label{pf.lem.det.laplace-multi.Apq.3b.pf.1c}%
\end{align}
\par
Recall that $\left(  p_{1},p_{2},\ldots,p_{k}\right)  $ is the list of all
elements of $P$ in increasing order (with no repetitions). Hence, Lemma
\ref{lem.sol.Ialbe.inclist} \textbf{(a)} (applied to $P$, $k$ and $\left(
p_{1},p_{2},\ldots,p_{k}\right)  $ instead of $S$, $s$ and $\left(
c_{1},c_{2},\ldots,c_{s}\right)  $) shows that the map $\left[  k\right]
\rightarrow P,\ h\mapsto p_{h}$ is well-defined and a bijection. In
particular, this map is a bijection. Hence, we can substitute $p_{h}$ for $i$
in the product $\prod_{i\in P}a_{i,\sigma\left(  i\right)  }$. We thus obtain%
\begin{align*}
\prod_{i\in P}a_{i,\sigma\left(  i\right)  }  &  =\prod_{h\in\left[  k\right]
}a_{p_{h},\sigma\left(  p_{h}\right)  }=\underbrace{\prod_{i\in\left[
k\right]  }}_{\substack{=\prod_{i\in\left\{  1,2,\ldots,k\right\}
}\\\text{(since }\left[  k\right]  =\left\{  1,2,\ldots,k\right\}  \text{)}%
}}a_{p_{i},\sigma\left(  p_{i}\right)  }\\
&  \ \ \ \ \ \ \ \ \ \ \left(  \text{here, we have renamed the index }h\text{
as }i\text{ in the product}\right) \\
&  =\prod_{i\in\left\{  1,2,\ldots,k\right\}  }\underbrace{a_{p_{i}%
,\sigma\left(  p_{i}\right)  }}_{\substack{=a_{\gamma\left(  i\right)
,\left(  \sigma\circ\gamma\right)  \left(  i\right)  }\\\text{(by
(\ref{pf.lem.det.laplace-multi.Apq.3b.pf.1c}))}}}\\
&  =\prod_{i\in\left\{  1,2,\ldots,k\right\}  }a_{\gamma\left(  i\right)
,\left(  \sigma\circ\gamma\right)  \left(  i\right)  }.
\end{align*}
This proves (\ref{pf.lem.det.laplace-multi.Apq.3b}).} and%
\begin{equation}
\prod_{i\in\left\{  k+1,k+2,\ldots,n\right\}  }b_{\gamma\left(  i\right)
,\left(  \sigma\circ\gamma\right)  \left(  i\right)  }=\prod_{i\in
\widetilde{P}}b_{i,\sigma\left(  i\right)  }
\label{pf.lem.det.laplace-multi.Apq.3c}%
\end{equation}
\footnote{\textit{Proof of (\ref{pf.lem.det.laplace-multi.Apq.3c}):} Let
$\sigma\in S_{n}$.
\par
We shall use the notation introduced in Definition \ref{def.sol.Ialbe.12n}.
Thus, $\left[  n-k\right]  =\left\{  1,2,\ldots,n-k\right\}  $ (by the
definition of $\left[  n-k\right]  $).
\par
Every $i\in\left\{  1,2,\ldots,n-k\right\}  $ satisfies%
\begin{equation}
r_{i}=\gamma\left(  k+i\right)  \ \ \ \ \ \ \ \ \ \ \left(  \text{by
(\ref{pf.lem.det.laplace-multi.Apq.gamma-tuple2})}\right)
\label{pf.lem.det.laplace-multi.Apq.3c.pf.1a}%
\end{equation}
and%
\begin{equation}
\sigma\left(  \underbrace{r_{i}}_{=\gamma\left(  k+i\right)  }\right)
=\sigma\left(  \gamma\left(  k+i\right)  \right)  =\left(  \sigma\circ
\gamma\right)  \left(  k+i\right)  .
\label{pf.lem.det.laplace-multi.Apq.3c.pf.1b}%
\end{equation}
Hence, every $i\in\left\{  1,2,\ldots,n-k\right\}  $ satisfies%
\begin{align}
b_{r_{i},\sigma\left(  r_{i}\right)  }  &  =b_{r_{i},\left(  \sigma\circ
\gamma\right)  \left(  k+i\right)  }\ \ \ \ \ \ \ \ \ \ \left(  \text{since
}\sigma\left(  r_{i}\right)  =\left(  \sigma\circ\gamma\right)  \left(
k+i\right)  \text{ (by (\ref{pf.lem.det.laplace-multi.Apq.3c.pf.1b}))}\right)
\nonumber\\
&  =b_{\gamma\left(  k+i\right)  ,\left(  \sigma\circ\gamma\right)  \left(
k+i\right)  }\ \ \ \ \ \ \ \ \ \ \left(  \text{since }r_{i}=\gamma\left(
k+i\right)  \text{ (by (\ref{pf.lem.det.laplace-multi.Apq.3c.pf.1a}))}\right)
\nonumber\\
&  =b_{\gamma\left(  i+k\right)  ,\left(  \sigma\circ\gamma\right)  \left(
i+k\right)  }\ \ \ \ \ \ \ \ \ \ \left(  \text{since }k+i=i+k\right)  .
\label{pf.lem.det.laplace-multi.Apq.3c.pf.1c}%
\end{align}
\par
Recall that $\left(  r_{1},r_{2},\ldots,r_{n-k}\right)  $ is the list of all
elements of $\widetilde{P}$ in increasing order (with no repetitions). Hence,
Lemma \ref{lem.sol.Ialbe.inclist} \textbf{(a)} (applied to $\widetilde{P}$,
$n-k$ and $\left(  r_{1},r_{2},\ldots,r_{n-k}\right)  $ instead of $S$, $s$
and $\left(  c_{1},c_{2},\ldots,c_{s}\right)  $) shows that the map $\left[
n-k\right]  \rightarrow\widetilde{P},\ h\mapsto r_{h}$ is well-defined and a
bijection. In particular, this map is a bijection. Hence, we can substitute
$r_{h}$ for $i$ in the product $\prod_{i\in\widetilde{P}}b_{i,\sigma\left(
i\right)  }$. We thus obtain%
\begin{align*}
\prod_{i\in\widetilde{P}}b_{i,\sigma\left(  i\right)  }  &  =\prod
_{h\in\left[  n-k\right]  }b_{r_{h},\sigma\left(  r_{h}\right)  }%
=\underbrace{\prod_{i\in\left[  n-k\right]  }}_{\substack{=\prod_{i\in\left\{
1,2,\ldots,n-k\right\}  }\\\text{(since }\left[  n-k\right]  =\left\{
1,2,\ldots,n-k\right\}  \text{)}}}b_{r_{i},\sigma\left(  r_{i}\right)  }\\
&  \ \ \ \ \ \ \ \ \ \ \left(  \text{here, we have renamed the index }h\text{
as }i\text{ in the product}\right) \\
&  =\underbrace{\prod_{i\in\left\{  1,2,\ldots,n-k\right\}  }}_{=\prod
_{i=1}^{n-k}}\underbrace{b_{r_{i},\sigma\left(  r_{i}\right)  }}%
_{\substack{=b_{\gamma\left(  i+k\right)  ,\left(  \sigma\circ\gamma\right)
\left(  i+k\right)  }\\\text{(by (\ref{pf.lem.det.laplace-multi.Apq.3c.pf.1c}%
))}}}\\
&  =\prod_{i=1}^{n-k}b_{\gamma\left(  i+k\right)  ,\left(  \sigma\circ
\gamma\right)  \left(  i+k\right)  }=\underbrace{\prod_{i=k+1}^{n}}%
_{=\prod_{i\in\left\{  k+1,k+2,\ldots,n\right\}  }}b_{\gamma\left(  i\right)
,\left(  \sigma\circ\gamma\right)  \left(  i\right)  }\\
&  \ \ \ \ \ \ \ \ \ \ \left(  \text{here, we have substituted }i\text{ for
}i+k\text{ in the product}\right) \\
&  =\prod_{i\in\left\{  k+1,k+2,\ldots,n\right\}  }b_{\gamma\left(  i\right)
,\left(  \sigma\circ\gamma\right)  \left(  i\right)  }.
\end{align*}
This proves (\ref{pf.lem.det.laplace-multi.Apq.3c}).}.
\end{verlong}

\begin{vershort}
Now, Lemma \ref{lem.sol.det.laplace-multi.1} (applied to $A^{\prime}$,
$a_{\gamma\left(  i\right)  ,j}$, $B^{\prime}$ and $b_{\gamma\left(  i\right)
,j}$ instead of $A$, $a_{i,j}$, $B$ and $b_{i,j}$) yields%
\begin{align*}
&  \sum_{\substack{\sigma\in S_{n};\\\sigma\left(  \left\{  1,2,\ldots
,k\right\}  \right)  =Q}}\left(  -1\right)  ^{\sigma}\left(  \prod
_{i\in\left\{  1,2,\ldots,k\right\}  }a_{\gamma\left(  i\right)
,\sigma\left(  i\right)  }\right)  \left(  \prod_{i\in\left\{  k+1,k+2,\ldots
,n\right\}  }b_{\gamma\left(  i\right)  ,\sigma\left(  i\right)  }\right) \\
&  =\left(  -1\right)  ^{\left(  1+2+\cdots+k\right)  +\sum Q}\det
\underbrace{\left(  \operatorname*{sub}\nolimits_{\left(  1,2,\ldots,k\right)
}^{w\left(  Q\right)  }\left(  A^{\prime}\right)  \right)  }%
_{\substack{=\operatorname*{sub}\nolimits_{w\left(  P\right)  }^{w\left(
Q\right)  }A\\\text{(by (\ref{pf.lem.det.laplace-multi.Apq.sub1}))}%
}}\underbrace{\det\left(  \operatorname*{sub}\nolimits_{\left(  k+1,k+2,\ldots
,n\right)  }^{w\left(  \widetilde{Q}\right)  }\left(  B^{\prime}\right)
\right)  }_{\substack{=\operatorname*{sub}\nolimits_{w\left(  \widetilde{P}%
\right)  }^{w\left(  \widetilde{Q}\right)  }B\\\text{(by
(\ref{pf.lem.det.laplace-multi.Apq.sub2}))}}}\\
&  =\left(  -1\right)  ^{\left(  1+2+\cdots+k\right)  +\sum Q}\det\left(
\operatorname*{sub}\nolimits_{w\left(  P\right)  }^{w\left(  Q\right)
}A\right)  \det\left(  \operatorname*{sub}\nolimits_{w\left(  \widetilde{P}%
\right)  }^{w\left(  \widetilde{Q}\right)  }B\right)  .
\end{align*}
Comparing this with%
\begin{align*}
&  \sum_{\substack{\sigma\in S_{n};\\\sigma\left(  \left\{  1,2,\ldots
,k\right\}  \right)  =Q}}\left(  -1\right)  ^{\sigma}\left(  \prod
_{i\in\left\{  1,2,\ldots,k\right\}  }a_{\gamma\left(  i\right)
,\sigma\left(  i\right)  }\right)  \left(  \prod_{i\in\left\{  k+1,k+2,\ldots
,n\right\}  }b_{\gamma\left(  i\right)  ,\sigma\left(  i\right)  }\right) \\
&  =\underbrace{\sum_{\substack{\sigma\in S_{n};\\\left(  \sigma\circ
\gamma\right)  \left(  \left\{  1,2,\ldots,k\right\}  \right)  =Q}%
}}_{\substack{=\sum_{\substack{\sigma\in S_{n};\\\sigma\left(  P\right)
=Q}}\\\text{(since every }\sigma\in S_{n}\\\text{satisfies }\left(
\sigma\circ\gamma\right)  \left(  \left\{  1,2,\ldots,k\right\}  \right)
=\sigma\left(  P\right)  \\\text{(by (\ref{pf.lem.det.laplace-multi.Apq.3a}%
)))}}}\underbrace{\left(  -1\right)  ^{\sigma\circ\gamma}}_{\substack{=\left(
-1\right)  ^{\sigma}\cdot\left(  -1\right)  ^{\gamma}\\\text{(by
(\ref{eq.sign.prod}), applied}\\\text{to }\tau=\gamma\text{)}}}\\
&  \ \ \ \ \ \ \ \ \ \ \underbrace{\left(  \prod_{i\in\left\{  1,2,\ldots
,k\right\}  }a_{\gamma\left(  i\right)  ,\left(  \sigma\circ\gamma\right)
\left(  i\right)  }\right)  }_{\substack{=\prod_{i\in P}a_{i,\sigma\left(
i\right)  }\\\text{(by (\ref{pf.lem.det.laplace-multi.Apq.short.3b}))}%
}}\underbrace{\left(  \prod_{i\in\left\{  k+1,k+2,\ldots,n\right\}  }%
b_{\gamma\left(  i\right)  ,\left(  \sigma\circ\gamma\right)  \left(
i\right)  }\right)  }_{\substack{=\prod_{i\in\widetilde{P}}b_{i,\sigma\left(
i\right)  }\\\text{(by (\ref{pf.lem.det.laplace-multi.Apq.short.3c}))}}}\\
&  \ \ \ \ \ \ \ \ \ \ \left(
\begin{array}
[c]{c}%
\text{here, we have substituted }\sigma\circ\gamma\text{ for }\sigma\text{ in
the sum, since}\\
\text{the map }S_{n}\rightarrow S_{n},\ \sigma\mapsto\sigma\circ\gamma\text{
is a bijection}%
\end{array}
\right) \\
&  =\sum_{\substack{\sigma\in S_{n};\\\sigma\left(  P\right)  =Q}}\left(
-1\right)  ^{\sigma}\cdot\left(  -1\right)  ^{\gamma}\left(  \prod_{i\in
P}a_{i,\sigma\left(  i\right)  }\right)  \left(  \prod_{i\in\widetilde{P}%
}b_{i,\sigma\left(  i\right)  }\right) \\
&  =\left(  -1\right)  ^{\gamma}\cdot\sum_{\substack{\sigma\in S_{n}%
;\\\sigma\left(  P\right)  =Q}}\left(  -1\right)  ^{\sigma}\left(  \prod_{i\in
P}a_{i,\sigma\left(  i\right)  }\right)  \left(  \prod_{i\in\widetilde{P}%
}b_{i,\sigma\left(  i\right)  }\right)  ,
\end{align*}
we obtain%
\begin{align*}
&  \left(  -1\right)  ^{\left(  1+2+\cdots+k\right)  +\sum Q}\det\left(
\operatorname*{sub}\nolimits_{w\left(  P\right)  }^{w\left(  Q\right)
}A\right)  \det\left(  \operatorname*{sub}\nolimits_{w\left(  \widetilde{P}%
\right)  }^{w\left(  \widetilde{Q}\right)  }B\right) \\
&  =\left(  -1\right)  ^{\gamma}\cdot\sum_{\substack{\sigma\in S_{n}%
;\\\sigma\left(  P\right)  =Q}}\left(  -1\right)  ^{\sigma}\left(  \prod_{i\in
P}a_{i,\sigma\left(  i\right)  }\right)  \left(  \prod_{i\in\widetilde{P}%
}b_{i,\sigma\left(  i\right)  }\right)  .
\end{align*}
Multiplying both sides of this equality by $\left(  -1\right)  ^{\sum
P-\left(  1+2+\cdots+k\right)  }$, we obtain%
\begin{align*}
&  \left(  -1\right)  ^{\sum P-\left(  1+2+\cdots+k\right)  }\left(
-1\right)  ^{\left(  1+2+\cdots+k\right)  +\sum Q}\det\left(
\operatorname*{sub}\nolimits_{w\left(  P\right)  }^{w\left(  Q\right)
}A\right)  \det\left(  \operatorname*{sub}\nolimits_{w\left(  \widetilde{P}%
\right)  }^{w\left(  \widetilde{Q}\right)  }B\right) \\
&  =\underbrace{\left(  -1\right)  ^{\sum P-\left(  1+2+\cdots+k\right)
}\left(  -1\right)  ^{\gamma}}_{\substack{=1\\\text{(by
(\ref{pf.lem.det.laplace-multi.Apq.sign-gamma-multed}))}}}\cdot\sum
_{\substack{\sigma\in S_{n};\\\sigma\left(  P\right)  =Q}}\left(  -1\right)
^{\sigma}\left(  \prod_{i\in P}a_{i,\sigma\left(  i\right)  }\right)  \left(
\prod_{i\in\widetilde{P}}b_{i,\sigma\left(  i\right)  }\right) \\
&  =\sum_{\substack{\sigma\in S_{n};\\\sigma\left(  P\right)  =Q}}\left(
-1\right)  ^{\sigma}\left(  \prod_{i\in P}a_{i,\sigma\left(  i\right)
}\right)  \left(  \prod_{i\in\widetilde{P}}b_{i,\sigma\left(  i\right)
}\right)  .
\end{align*}
Thus,%
\begin{align*}
&  \sum_{\substack{\sigma\in S_{n};\\\sigma\left(  P\right)  =Q}}\left(
-1\right)  ^{\sigma}\left(  \prod_{i\in P}a_{i,\sigma\left(  i\right)
}\right)  \left(  \prod_{i\in\widetilde{P}}b_{i,\sigma\left(  i\right)
}\right) \\
&  =\underbrace{\left(  -1\right)  ^{\sum P-\left(  1+2+\cdots+k\right)
}\left(  -1\right)  ^{\left(  1+2+\cdots+k\right)  +\sum Q}}%
_{\substack{=\left(  -1\right)  ^{\left(  \sum P-\left(  1+2+\cdots+k\right)
\right)  +\left(  \left(  1+2+\cdots+k\right)  +\sum Q\right)  }\\=\left(
-1\right)  ^{\sum P+\sum Q}\\\text{(since }\left(  \sum P-\left(
1+2+\cdots+k\right)  \right)  +\left(  \left(  1+2+\cdots+k\right)  +\sum
Q\right)  =\sum P+\sum Q\text{)}}}\det\left(  \operatorname*{sub}%
\nolimits_{w\left(  P\right)  }^{w\left(  Q\right)  }A\right)  \det\left(
\operatorname*{sub}\nolimits_{w\left(  \widetilde{P}\right)  }^{w\left(
\widetilde{Q}\right)  }B\right) \\
&  =\left(  -1\right)  ^{\sum P+\sum Q}\det\left(  \operatorname*{sub}%
\nolimits_{w\left(  P\right)  }^{w\left(  Q\right)  }A\right)  \det\left(
\operatorname*{sub}\nolimits_{w\left(  \widetilde{P}\right)  }^{w\left(
\widetilde{Q}\right)  }B\right)  .
\end{align*}
This proves Lemma \ref{lem.det.laplace-multi.Apq}. \qedhere

\end{vershort}

\begin{verlong}
Now, Lemma \ref{lem.sol.det.laplace-multi.1} (applied to $A^{\prime}$,
$a_{\gamma\left(  i\right)  ,j}$, $B^{\prime}$ and $b_{\gamma\left(  i\right)
,j}$ instead of $A$, $a_{i,j}$, $B$ and $b_{i,j}$) yields%
\begin{align*}
&  \sum_{\substack{\sigma\in S_{n};\\\sigma\left(  \left\{  1,2,\ldots
,k\right\}  \right)  =Q}}\left(  -1\right)  ^{\sigma}\left(  \prod
_{i\in\left\{  1,2,\ldots,k\right\}  }a_{\gamma\left(  i\right)
,\sigma\left(  i\right)  }\right)  \left(  \prod_{i\in\left\{  k+1,k+2,\ldots
,n\right\}  }b_{\gamma\left(  i\right)  ,\sigma\left(  i\right)  }\right) \\
&  =\left(  -1\right)  ^{\left(  1+2+\cdots+k\right)  +\sum Q}\det
\underbrace{\left(  \operatorname*{sub}\nolimits_{\left(  1,2,\ldots,k\right)
}^{w\left(  Q\right)  }\left(  A^{\prime}\right)  \right)  }%
_{\substack{=\operatorname*{sub}\nolimits_{w\left(  P\right)  }^{w\left(
Q\right)  }A\\\text{(by (\ref{pf.lem.det.laplace-multi.Apq.sub1}))}%
}}\underbrace{\det\left(  \operatorname*{sub}\nolimits_{\left(  k+1,k+2,\ldots
,n\right)  }^{w\left(  \widetilde{Q}\right)  }\left(  B^{\prime}\right)
\right)  }_{\substack{=\operatorname*{sub}\nolimits_{w\left(  \widetilde{P}%
\right)  }^{w\left(  \widetilde{Q}\right)  }B\\\text{(by
(\ref{pf.lem.det.laplace-multi.Apq.sub2}))}}}\\
&  =\left(  -1\right)  ^{\left(  1+2+\cdots+k\right)  +\sum Q}\det\left(
\operatorname*{sub}\nolimits_{w\left(  P\right)  }^{w\left(  Q\right)
}A\right)  \det\left(  \operatorname*{sub}\nolimits_{w\left(  \widetilde{P}%
\right)  }^{w\left(  \widetilde{Q}\right)  }B\right)  .
\end{align*}
Comparing this with%
\begin{align*}
&  \sum_{\substack{\sigma\in S_{n};\\\sigma\left(  \left\{  1,2,\ldots
,k\right\}  \right)  =Q}}\left(  -1\right)  ^{\sigma}\left(  \prod
_{i\in\left\{  1,2,\ldots,k\right\}  }a_{\gamma\left(  i\right)
,\sigma\left(  i\right)  }\right)  \left(  \prod_{i\in\left\{  k+1,k+2,\ldots
,n\right\}  }b_{\gamma\left(  i\right)  ,\sigma\left(  i\right)  }\right) \\
&  =\underbrace{\sum_{\substack{\sigma\in S_{n};\\\left(  \sigma\circ
\gamma\right)  \left(  \left\{  1,2,\ldots,k\right\}  \right)  =Q}%
}}_{\substack{=\sum_{\substack{\sigma\in S_{n};\\\sigma\left(  P\right)
=Q}}\\\text{(since every }\sigma\in S_{n}\\\text{satisfies }\left(
\sigma\circ\gamma\right)  \left(  \left\{  1,2,\ldots,k\right\}  \right)
=\sigma\left(  P\right)  \\\text{(by (\ref{pf.lem.det.laplace-multi.Apq.3a}%
)))}}}\underbrace{\left(  -1\right)  ^{\sigma\circ\gamma}}_{\substack{=\left(
-1\right)  ^{\sigma}\cdot\left(  -1\right)  ^{\gamma}\\\text{(by
(\ref{eq.sign.prod}), applied}\\\text{to }\tau=\gamma\text{)}}}\\
&  \ \ \ \ \ \ \ \ \ \ \underbrace{\left(  \prod_{i\in\left\{  1,2,\ldots
,k\right\}  }a_{\gamma\left(  i\right)  ,\left(  \sigma\circ\gamma\right)
\left(  i\right)  }\right)  }_{\substack{=\prod_{i\in P}a_{i,\sigma\left(
i\right)  }\\\text{(by (\ref{pf.lem.det.laplace-multi.Apq.3b}))}%
}}\underbrace{\left(  \prod_{i\in\left\{  k+1,k+2,\ldots,n\right\}  }%
b_{\gamma\left(  i\right)  ,\left(  \sigma\circ\gamma\right)  \left(
i\right)  }\right)  }_{\substack{=\prod_{i\in\widetilde{P}}b_{i,\sigma\left(
i\right)  }\\\text{(by (\ref{pf.lem.det.laplace-multi.Apq.3c}))}}}\\
&  \ \ \ \ \ \ \ \ \ \ \left(
\begin{array}
[c]{c}%
\text{here, we have substituted }\sigma\circ\gamma\text{ for }\sigma\text{ in
the sum, since}\\
\text{the map }S_{n}\rightarrow S_{n},\ \sigma\mapsto\sigma\circ\gamma\text{
is a bijection}%
\end{array}
\right) \\
&  =\sum_{\substack{\sigma\in S_{n};\\\sigma\left(  P\right)  =Q}}\left(
-1\right)  ^{\sigma}\cdot\left(  -1\right)  ^{\gamma}\left(  \prod_{i\in
P}a_{i,\sigma\left(  i\right)  }\right)  \left(  \prod_{i\in\widetilde{P}%
}b_{i,\sigma\left(  i\right)  }\right) \\
&  =\left(  -1\right)  ^{\gamma}\cdot\sum_{\substack{\sigma\in S_{n}%
;\\\sigma\left(  P\right)  =Q}}\left(  -1\right)  ^{\sigma}\left(  \prod_{i\in
P}a_{i,\sigma\left(  i\right)  }\right)  \left(  \prod_{i\in\widetilde{P}%
}b_{i,\sigma\left(  i\right)  }\right)  ,
\end{align*}
we obtain%
\begin{align*}
&  \left(  -1\right)  ^{\left(  1+2+\cdots+k\right)  +\sum Q}\det\left(
\operatorname*{sub}\nolimits_{w\left(  P\right)  }^{w\left(  Q\right)
}A\right)  \det\left(  \operatorname*{sub}\nolimits_{w\left(  \widetilde{P}%
\right)  }^{w\left(  \widetilde{Q}\right)  }B\right) \\
&  =\left(  -1\right)  ^{\gamma}\cdot\sum_{\substack{\sigma\in S_{n}%
;\\\sigma\left(  P\right)  =Q}}\left(  -1\right)  ^{\sigma}\left(  \prod_{i\in
P}a_{i,\sigma\left(  i\right)  }\right)  \left(  \prod_{i\in\widetilde{P}%
}b_{i,\sigma\left(  i\right)  }\right)  .
\end{align*}
Multiplying both sides of this equality by $\left(  -1\right)  ^{\sum
P-\left(  1+2+\cdots+k\right)  }$, we obtain%
\begin{align*}
&  \left(  -1\right)  ^{\sum P-\left(  1+2+\cdots+k\right)  }\left(
-1\right)  ^{\left(  1+2+\cdots+k\right)  +\sum Q}\det\left(
\operatorname*{sub}\nolimits_{w\left(  P\right)  }^{w\left(  Q\right)
}A\right)  \det\left(  \operatorname*{sub}\nolimits_{w\left(  \widetilde{P}%
\right)  }^{w\left(  \widetilde{Q}\right)  }B\right) \\
&  =\underbrace{\left(  -1\right)  ^{\sum P-\left(  1+2+\cdots+k\right)
}\left(  -1\right)  ^{\gamma}}_{\substack{=1\\\text{(by
(\ref{pf.lem.det.laplace-multi.Apq.sign-gamma-multed}))}}}\cdot\sum
_{\substack{\sigma\in S_{n};\\\sigma\left(  P\right)  =Q}}\left(  -1\right)
^{\sigma}\left(  \prod_{i\in P}a_{i,\sigma\left(  i\right)  }\right)  \left(
\prod_{i\in\widetilde{P}}b_{i,\sigma\left(  i\right)  }\right) \\
&  =\sum_{\substack{\sigma\in S_{n};\\\sigma\left(  P\right)  =Q}}\left(
-1\right)  ^{\sigma}\left(  \prod_{i\in P}a_{i,\sigma\left(  i\right)
}\right)  \left(  \prod_{i\in\widetilde{P}}b_{i,\sigma\left(  i\right)
}\right)  .
\end{align*}
Thus,%
\begin{align*}
&  \sum_{\substack{\sigma\in S_{n};\\\sigma\left(  P\right)  =Q}}\left(
-1\right)  ^{\sigma}\left(  \prod_{i\in P}a_{i,\sigma\left(  i\right)
}\right)  \left(  \prod_{i\in\widetilde{P}}b_{i,\sigma\left(  i\right)
}\right) \\
&  =\underbrace{\left(  -1\right)  ^{\sum P-\left(  1+2+\cdots+k\right)
}\left(  -1\right)  ^{\left(  1+2+\cdots+k\right)  +\sum Q}}%
_{\substack{=\left(  -1\right)  ^{\left(  \sum P-\left(  1+2+\cdots+k\right)
\right)  +\left(  \left(  1+2+\cdots+k\right)  +\sum Q\right)  }\\=\left(
-1\right)  ^{\sum P+\sum Q}\\\text{(since }\left(  \sum P-\left(
1+2+\cdots+k\right)  \right)  +\left(  \left(  1+2+\cdots+k\right)  +\sum
Q\right)  =\sum P+\sum Q\text{)}}}\det\left(  \operatorname*{sub}%
\nolimits_{w\left(  P\right)  }^{w\left(  Q\right)  }A\right)  \det\left(
\operatorname*{sub}\nolimits_{w\left(  \widetilde{P}\right)  }^{w\left(
\widetilde{Q}\right)  }B\right) \\
&  =\left(  -1\right)  ^{\sum P+\sum Q}\det\left(  \operatorname*{sub}%
\nolimits_{w\left(  P\right)  }^{w\left(  Q\right)  }A\right)  \det\left(
\operatorname*{sub}\nolimits_{w\left(  \widetilde{P}\right)  }^{w\left(
\widetilde{Q}\right)  }B\right)  .
\end{align*}
This proves Lemma \ref{lem.det.laplace-multi.Apq}.
\end{verlong}
\end{proof}

For the sake of convenience, let us restate a simplified particular case of
Lemma \ref{lem.det.laplace-multi.Apq} for $A=B$:

\begin{lemma}
\label{lem.sol.det.laplace-multi.2}Let $n\in\mathbb{N}$. For any subset $I$ of
$\left\{  1,2,\ldots,n\right\}  $, we let $\widetilde{I}$ denote the
complement $\left\{  1,2,\ldots,n\right\}  \setminus I$ of $I$.

Let $A=\left(  a_{i,j}\right)  _{1\leq i\leq n,\ 1\leq j\leq n}$ be an
$n\times n$-matrix. Let $P$ and $Q$ be two subsets of $\left\{  1,2,\ldots
,n\right\}  $ such that $\left\vert P\right\vert =\left\vert Q\right\vert $.
Then,%
\[
\sum_{\substack{\sigma\in S_{n};\\\sigma\left(  P\right)  =Q}}\left(
-1\right)  ^{\sigma}\prod_{i=1}^{n}a_{i,\sigma\left(  i\right)  }=\left(
-1\right)  ^{\sum P+\sum Q}\det\left(  \operatorname*{sub}\nolimits_{w\left(
P\right)  }^{w\left(  Q\right)  }A\right)  \det\left(  \operatorname*{sub}%
\nolimits_{w\left(  \widetilde{P}\right)  }^{w\left(  \widetilde{Q}\right)
}A\right)  .
\]

\end{lemma}

\begin{proof}
[Proof of Lemma \ref{lem.sol.det.laplace-multi.2}.]Every $\sigma\in S_{n}$
satisfies%
\begin{equation}
\prod_{i=1}^{n}a_{i,\sigma\left(  i\right)  }=\left(  \prod_{i\in
P}a_{i,\sigma\left(  i\right)  }\right)  \left(  \prod_{i\in\widetilde{P}%
}a_{i,\sigma\left(  i\right)  }\right)
\label{pf.lem.sol.det.laplace-multi.2.1}%
\end{equation}
\footnote{\textit{Proof of (\ref{pf.lem.sol.det.laplace-multi.2.1}):} Let
$\sigma\in S_{n}$. Notice that $\widetilde{P}=\left\{  1,2,\ldots,n\right\}
\setminus P$ (by the definition of $\widetilde{P}$). Now,%
\begin{align*}
\underbrace{\prod_{i=1}^{n}}_{=\prod_{i\in\left\{  1,2,\ldots,n\right\}  }%
}a_{i,\sigma\left(  i\right)  }  &  =\prod_{i\in\left\{  1,2,\ldots,n\right\}
}a_{i,\sigma\left(  i\right)  }=\left(  \underbrace{\prod_{\substack{i\in
\left\{  1,2,\ldots,n\right\}  ;\\i\in P}}}_{\substack{=\prod_{i\in
P}\\\text{(since }P\subseteq\left\{  1,2,\ldots,n\right\}  \text{)}%
}}a_{i,\sigma\left(  i\right)  }\right)  \left(  \underbrace{\prod
_{\substack{i\in\left\{  1,2,\ldots,n\right\}  ;\\i\notin P}}}%
_{\substack{=\prod_{i\in\left\{  1,2,\ldots,n\right\}  \setminus P}%
=\prod_{i\in\widetilde{P}}\\\text{(since }\left\{  1,2,\ldots,n\right\}
\setminus P=\widetilde{P}\text{)}}}a_{i,\sigma\left(  i\right)  }\right) \\
&  \ \ \ \ \ \ \ \ \ \ \left(
\begin{array}
[c]{c}%
\text{since every }i\in\left\{  1,2,\ldots,n\right\}  \text{ satisfies either
}i\in P\\
\text{or }i\notin P\text{ (but not both)}%
\end{array}
\right) \\
&  =\left(  \prod_{i\in P}a_{i,\sigma\left(  i\right)  }\right)  \left(
\prod_{i\in\widetilde{P}}a_{i,\sigma\left(  i\right)  }\right)  .
\end{align*}
This proves (\ref{pf.lem.sol.det.laplace-multi.2.1}).}. Now,%
\begin{align*}
&  \sum_{\substack{\sigma\in S_{n};\\\sigma\left(  P\right)  =Q}}\left(
-1\right)  ^{\sigma}\underbrace{\prod_{i=1}^{n}a_{i,\sigma\left(  i\right)  }%
}_{\substack{=\left(  \prod_{i\in P}a_{i,\sigma\left(  i\right)  }\right)
\left(  \prod_{i\in\widetilde{P}}a_{i,\sigma\left(  i\right)  }\right)
\\\text{(by (\ref{pf.lem.sol.det.laplace-multi.2.1}))}}}\\
&  =\sum_{\substack{\sigma\in S_{n};\\\sigma\left(  P\right)  =Q}}\left(
-1\right)  ^{\sigma}\left(  \prod_{i\in P}a_{i,\sigma\left(  i\right)
}\right)  \left(  \prod_{i\in\widetilde{P}}a_{i,\sigma\left(  i\right)
}\right) \\
&  =\left(  -1\right)  ^{\sum P+\sum Q}\det\left(  \operatorname*{sub}%
\nolimits_{w\left(  P\right)  }^{w\left(  Q\right)  }A\right)  \det\left(
\operatorname*{sub}\nolimits_{w\left(  \widetilde{P}\right)  }^{w\left(
\widetilde{Q}\right)  }A\right)
\end{align*}
(by Lemma \ref{lem.det.laplace-multi.Apq} (applied to $B=A$ and $b_{i,j}%
=a_{i,j}$)). This proves Lemma \ref{lem.sol.det.laplace-multi.2}.
\end{proof}

Another fact that we will need is very simple:

\begin{lemma}
\label{lem.sol.det.laplace-multi.4}Let $n\in\mathbb{N}$. Let $\sigma\in S_{n}%
$. Let $P$ be a subset of $\left\{  1,2,\ldots,n\right\}  $.

\textbf{(a)} The set $\sigma\left(  P\right)  $ is a subset of $\left\{
1,2,\ldots,n\right\}  $ satisfying $\left\vert \sigma\left(  P\right)
\right\vert =\left\vert P\right\vert $.

\textbf{(b)} Let $Q$ be a subset of $\left\{  1,2,\ldots,n\right\}  $. Then,
$\sigma^{-1}\left(  Q\right)  =P$ holds if and only if $\sigma\left(
P\right)  =Q$.
\end{lemma}

\begin{vershort}
\begin{proof}
[Proof of Lemma \ref{lem.sol.det.laplace-multi.4}.]We have $\sigma\in S_{n}$.
In other words, $\sigma$ is a permutation of the set $\left\{  1,2,\ldots
,n\right\}  $. In other words, $\sigma$ is a bijective map $\left\{
1,2,\ldots,n\right\}  \rightarrow\left\{  1,2,\ldots,n\right\}  $. The map
$\sigma$ is bijective, and therefore injective.

\textbf{(a)} Lemma \ref{lem.jectivity.pigeon0} \textbf{(c)} (applied to
$U=\left\{  1,2,\ldots,n\right\}  $, $V=\left\{  1,2,\ldots,n\right\}  $,
$f=\sigma$ and $S=P$) shows that $\left\vert \sigma\left(  P\right)
\right\vert =\left\vert P\right\vert $ (since the map $\sigma$ is injective).
This proves Lemma \ref{lem.sol.det.laplace-multi.4} \textbf{(a)}.

\textbf{(b)} We have $\sigma\left(  \sigma^{-1}\left(  Q\right)  \right)  =Q$
(since $\sigma$ is bijective). But if $U$ and $V$ are two subsets of $\left\{
1,2,\ldots,n\right\}  $, then $U=V$ holds if and only if $\sigma\left(
U\right)  =\sigma\left(  V\right)  $ (because $\sigma$ is bijective). Applying
this to $U=\sigma^{-1}\left(  Q\right)  $ and $V=P$, we conclude that
$\sigma^{-1}\left(  Q\right)  =P$ holds if and only if $\sigma\left(
\sigma^{-1}\left(  Q\right)  \right)  =\sigma\left(  P\right)  $. Hence, we
have the following chain of equivalences:%
\[
\left(  \sigma^{-1}\left(  Q\right)  =P\right)  \ \Longleftrightarrow\ \left(
\underbrace{\sigma\left(  \sigma^{-1}\left(  Q\right)  \right)  }_{=Q}%
=\sigma\left(  P\right)  \right)  \ \Longleftrightarrow\ \left(
Q=\sigma\left(  P\right)  \right)  \ \Longleftrightarrow\ \left(
\sigma\left(  P\right)  =Q\right)  .
\]
In other words, $\sigma^{-1}\left(  Q\right)  =P$ holds if and only if
$\sigma\left(  P\right)  =Q$. Lemma \ref{lem.sol.det.laplace-multi.4}
\textbf{(b)} is now proven.
\end{proof}
\end{vershort}

\begin{verlong}
\begin{proof}
[Proof of Lemma \ref{lem.sol.det.laplace-multi.4}.]We have $\sigma\in S_{n}$.
In other words, $\sigma$ is a permutation of the set $\left\{  1,2,\ldots
,n\right\}  $ (since $S_{n}$ is the set of all permutations of the set
$\left\{  1,2,\ldots,n\right\}  $). In other words, $\sigma$ is a bijective
map $\left\{  1,2,\ldots,n\right\}  \rightarrow\left\{  1,2,\ldots,n\right\}
$. The map $\sigma$ is bijective, and therefore injective.

\textbf{(a)} Clearly, the set $\sigma\left(  P\right)  $ is a subset of
$\left\{  1,2,\ldots,n\right\}  $. It remains to prove that $\left\vert
\sigma\left(  P\right)  \right\vert =\left\vert P\right\vert $.

Lemma \ref{lem.jectivity.pigeon0} \textbf{(c)} (applied to $U=\left\{
1,2,\ldots,n\right\}  $, $V=\left\{  1,2,\ldots,n\right\}  $, $f=\sigma$ and
$S=P$) shows that $\left\vert \sigma\left(  P\right)  \right\vert =\left\vert
P\right\vert $ (since the map $\sigma$ is injective). This proves Lemma
\ref{lem.sol.det.laplace-multi.4} \textbf{(a)}.

\textbf{(b)} We shall first prove the logical implication%
\begin{equation}
\left(  \sigma^{-1}\left(  Q\right)  =P\right)  \ \Longrightarrow\ \left(
\sigma\left(  P\right)  =Q\right)  .
\label{pf.lem.sol.det.laplace-multi.4.dir1}%
\end{equation}

[\textit{Proof of (\ref{pf.lem.sol.det.laplace-multi.4.dir1}):} Assume that
$\sigma^{-1}\left(  Q\right)  =P$. We shall show that $\sigma\left(  P\right)
=Q$.

We have $\sigma\left(  P\right)  \subseteq Q$\ \ \ \ \footnote{\textit{Proof.}
Let $z\in\sigma\left(  P\right)  $. Thus, $z=\sigma\left(  p\right)  $ for
some $p\in P$. Consider this $p$. We have $p\in P=\sigma^{-1}\left(  Q\right)
$, so that $\sigma\left(  p\right)  \in Q$. Thus, $z=\sigma\left(  p\right)
\in Q$.
\par
Now, forget that we fixed $z$. We thus have shown that $z\in Q$ for every
$z\in\sigma\left(  P\right)  $. In other words, $\sigma\left(  P\right)
\subseteq Q$. Qed.} and $Q\subseteq\sigma\left(  P\right)  $%
\ \ \ \ \footnote{\textit{Proof.} Let $q\in Q$. Now, $\sigma^{-1}\left(
q\right)  $ is a well-defined element of $\left\{  1,2,\ldots,n\right\}  $
(since the map $\sigma$ is bijective). This element $\sigma^{-1}\left(
q\right)  $ belongs to $\sigma^{-1}\left(  Q\right)  $ (since $\sigma\left(
\sigma^{-1}\left(  q\right)  \right)  =q\in Q$). Thus, $\sigma^{-1}\left(
q\right)  \in\sigma^{-1}\left(  Q\right)  =P$. Now, $q=\sigma\left(
\underbrace{\sigma^{-1}\left(  q\right)  }_{\in P}\right)  \in\sigma\left(
P\right)  $.
\par
Now, forget that we fixed $q$. We thus have shown that $q\in\sigma\left(
P\right)  $ for every $q\in Q$. In other words, $Q\subseteq\sigma\left(
P\right)  $. Qed.}. Combining these two relations, we obtain $\sigma\left(
P\right)  =Q$.

Now, forget that we assumed that $\sigma^{-1}\left(  Q\right)  =P$. We thus
have proven that $\sigma\left(  P\right)  =Q$ under the assumption that
$\sigma^{-1}\left(  Q\right)  =P$. In other words, we have proven the
implication (\ref{pf.lem.sol.det.laplace-multi.4.dir1}).]

Next, we shall prove the logical implication%
\begin{equation}
\left(  \sigma\left(  P\right)  =Q\right)  \ \Longrightarrow\ \left(
\sigma^{-1}\left(  Q\right)  =P\right)  .
\label{pf.lem.sol.det.laplace-multi.4.dir2}%
\end{equation}

[\textit{Proof of (\ref{pf.lem.sol.det.laplace-multi.4.dir2}):} Assume that
$\sigma\left(  P\right)  =Q$. We shall show that $\sigma^{-1}\left(  Q\right)
=P$.

We have $\sigma^{-1}\left(  Q\right)  \subseteq P$%
\ \ \ \ \footnote{\textit{Proof.} Let $z\in\sigma^{-1}\left(  Q\right)  $.
Thus, $\sigma\left(  z\right)  \in Q=\sigma\left(  P\right)  $. In other
words, $\sigma\left(  z\right)  =\sigma\left(  p\right)  $ for some $p\in P$.
Consider this $p$.
\par
The map $\sigma$ is injective. Thus, from $\sigma\left(  z\right)
=\sigma\left(  p\right)  $, we obtain $z=p$. Hence, $z=p\in P$.
\par
Now, forget that we fixed $z$. We thus have shown that $z\in P$ for every
$z\in\sigma^{-1}\left(  Q\right)  $. In other words, $\sigma^{-1}\left(
Q\right)  \subseteq P$. Qed.} and $P\subseteq\sigma^{-1}\left(  Q\right)
$\ \ \ \ \footnote{\textit{Proof.} Let $p\in P$. Then, $\sigma\left(
\underbrace{p}_{\in P}\right)  \in\sigma\left(  P\right)  =Q$, so that
$p\in\sigma^{-1}\left(  Q\right)  $.
\par
Now, forget that we fixed $p$. We thus have shown that $p\in\sigma^{-1}\left(
Q\right)  $ for every $p\in P$. In other words, $P\subseteq\sigma^{-1}\left(
Q\right)  $. Qed.}. Combining these two relations, we obtain $\sigma
^{-1}\left(  Q\right)  =P$.

Now, forget that we assumed that $\sigma\left(  P\right)  =Q$. We thus have
proven that $\left(  \sigma^{-1}\left(  Q\right)  =P\right)  $ under the
assumption that $\sigma\left(  P\right)  =Q$. In other words, we have proven
the implication (\ref{pf.lem.sol.det.laplace-multi.4.dir2}).]

Combining the two implications (\ref{pf.lem.sol.det.laplace-multi.4.dir1}) and
(\ref{pf.lem.sol.det.laplace-multi.4.dir2}), we obtain the logical equivalence
$\left(  \sigma^{-1}\left(  Q\right)  =P\right)  \ \Longleftrightarrow
\ \left(  \sigma\left(  P\right)  =Q\right)  $. In other words, $\sigma
^{-1}\left(  Q\right)  =P$ holds if and only if $\sigma\left(  P\right)  =Q$.
This proves Lemma \ref{lem.sol.det.laplace-multi.4} \textbf{(b)}.
\end{proof}
\end{verlong}

We can now step to the proof of Theorem \ref{thm.det.laplace-multi}:

\begin{proof}
[Proof of Theorem \ref{thm.det.laplace-multi}.]Write the $n\times n$-matrix
$A$ in the form $A=\left(  a_{i,j}\right)  _{1\leq i\leq n,\ 1\leq j\leq n}$.

If $P$ and $Q$ are two subsets of $\left\{  1,2,\ldots,n\right\}  $ satisfying
$\left\vert Q\right\vert \neq\left\vert P\right\vert $, then%
\begin{equation}
\sum_{\substack{\sigma\in S_{n};\\\sigma\left(  P\right)  =Q}}\left(
-1\right)  ^{\sigma}\prod_{i=1}^{n}a_{i,\sigma\left(  i\right)  }=0
\label{pf.thm.det.laplace-multi.emptysum}%
\end{equation}
\footnote{\textit{Proof of (\ref{pf.thm.det.laplace-multi.emptysum}):} Let $P$
and $Q$ be two subsets of $\left\{  1,2,\ldots,n\right\}  $ satisfying
$\left\vert Q\right\vert \neq\left\vert P\right\vert $.
\par
Let $\sigma\in S_{n}$ be such that $\sigma\left(  P\right)  =Q$. We shall
derive a contradiction.
\par
Indeed, Lemma \ref{lem.sol.det.laplace-multi.4} \textbf{(a)} shows that the
set $\sigma\left(  P\right)  $ is a subset of $\left\{  1,2,\ldots,n\right\}
$ satisfying $\left\vert \sigma\left(  P\right)  \right\vert =\left\vert
P\right\vert $. Hence, $\left\vert P\right\vert =\left\vert \underbrace{\sigma
\left(  P\right)  }_{=Q}\right\vert =\left\vert Q\right\vert \neq\left\vert
P\right\vert $. This is absurd. Hence, we have found a contradiction.
\par
Now, forget that we fixed $\sigma$. We thus have found a contradiction for
every $\sigma\in S_{n}$ satisfying $\sigma\left(  P\right)  =Q$. Thus, there
exists no $\sigma\in S_{n}$ satisfying $\sigma\left(  P\right)  =Q$. Hence,
the sum $\sum_{\substack{\sigma\in S_{n};\\\sigma\left(  P\right)  =Q}}\left(
-1\right)  ^{\sigma}\prod_{i=1}^{n}a_{i,\sigma\left(  i\right)  }$ is an empty
sum. Thus,%
\[
\sum_{\substack{\sigma\in S_{n};\\\sigma\left(  P\right)  =Q}}\left(
-1\right)  ^{\sigma}\prod_{i=1}^{n}a_{i,\sigma\left(  i\right)  }=\left(
\text{empty sum}\right)  =0.
\]
This proves (\ref{pf.thm.det.laplace-multi.emptysum}).}.

\textbf{(a)} Let $P$ be a subset of $\left\{  1,2,\ldots,n\right\}  $. Then,
(\ref{eq.det.eq.2}) yields%
\begin{align*}
\det A  &  =\underbrace{\sum_{\sigma\in S_{n}}}_{\substack{=\sum
_{Q\subseteq\left\{  1,2,\ldots,n\right\}  }\sum_{\substack{\sigma\in
S_{n};\\\sigma\left(  P\right)  =Q}}\\\text{(because for every }\sigma\in
S_{n}\text{, the set}\\\sigma\left(  P\right)  \text{ is a subset of }\left\{
1,2,\ldots,n\right\}  \text{)}}}\left(  -1\right)  ^{\sigma}\prod_{i=1}%
^{n}a_{i,\sigma\left(  i\right)  }=\sum_{Q\subseteq\left\{  1,2,\ldots
,n\right\}  }\sum_{\substack{\sigma\in S_{n};\\\sigma\left(  P\right)
=Q}}\left(  -1\right)  ^{\sigma}\prod_{i=1}^{n}a_{i,\sigma\left(  i\right)
}\\
&  =\sum_{\substack{Q\subseteq\left\{  1,2,\ldots,n\right\}  ;\\\left\vert
Q\right\vert =\left\vert P\right\vert }}\sum_{\substack{\sigma\in
S_{n};\\\sigma\left(  P\right)  =Q}}\left(  -1\right)  ^{\sigma}\prod
_{i=1}^{n}a_{i,\sigma\left(  i\right)  }+\sum_{\substack{Q\subseteq\left\{
1,2,\ldots,n\right\}  ;\\\left\vert Q\right\vert \neq\left\vert P\right\vert
}}\underbrace{\sum_{\substack{\sigma\in S_{n};\\\sigma\left(  P\right)
=Q}}\left(  -1\right)  ^{\sigma}\prod_{i=1}^{n}a_{i,\sigma\left(  i\right)  }%
}_{\substack{=0\\\text{(by (\ref{pf.thm.det.laplace-multi.emptysum}))}}}\\
&  \ \ \ \ \ \ \ \ \ \ \left(
\begin{array}
[c]{c}%
\text{since every subset }Q\text{ of }\left\{  1,2,\ldots,n\right\}  \text{
satisfies}\\
\text{either }\left\vert Q\right\vert =\left\vert P\right\vert \text{ or
}\left\vert Q\right\vert \neq\left\vert P\right\vert \text{ (but not both)}%
\end{array}
\right) \\
&  =\sum_{\substack{Q\subseteq\left\{  1,2,\ldots,n\right\}  ;\\\left\vert
Q\right\vert =\left\vert P\right\vert }}\sum_{\substack{\sigma\in
S_{n};\\\sigma\left(  P\right)  =Q}}\left(  -1\right)  ^{\sigma}\prod
_{i=1}^{n}a_{i,\sigma\left(  i\right)  }+\underbrace{\sum
_{\substack{Q\subseteq\left\{  1,2,\ldots,n\right\}  ;\\\left\vert
Q\right\vert \neq\left\vert P\right\vert }}0}_{=0}\\
&  =\sum_{\substack{Q\subseteq\left\{  1,2,\ldots,n\right\}  ;\\\left\vert
Q\right\vert =\left\vert P\right\vert }}\underbrace{\sum_{\substack{\sigma\in
S_{n};\\\sigma\left(  P\right)  =Q}}\left(  -1\right)  ^{\sigma}\prod
_{i=1}^{n}a_{i,\sigma\left(  i\right)  }}_{\substack{=\left(  -1\right)
^{\sum P+\sum Q}\det\left(  \operatorname*{sub}\nolimits_{w\left(  P\right)
}^{w\left(  Q\right)  }A\right)  \det\left(  \operatorname*{sub}%
\nolimits_{w\left(  \widetilde{P}\right)  }^{w\left(  \widetilde{Q}\right)
}A\right)  \\\text{(by Lemma \ref{lem.sol.det.laplace-multi.2}}\\\text{(since
}\left\vert P\right\vert =\left\vert Q\right\vert \text{ (since }\left\vert
Q\right\vert =\left\vert P\right\vert \text{)))}}}\\
&  =\sum_{\substack{Q\subseteq\left\{  1,2,\ldots,n\right\}  ;\\\left\vert
Q\right\vert =\left\vert P\right\vert }}\left(  -1\right)  ^{\sum P+\sum
Q}\det\left(  \operatorname*{sub}\nolimits_{w\left(  P\right)  }^{w\left(
Q\right)  }A\right)  \det\left(  \operatorname*{sub}\nolimits_{w\left(
\widetilde{P}\right)  }^{w\left(  \widetilde{Q}\right)  }A\right)  .
\end{align*}
This proves Theorem \ref{thm.det.laplace-multi} \textbf{(a)}.

\textbf{(b)} Let $Q$ be a subset of $\left\{  1,2,\ldots,n\right\}  $. Then,
(\ref{eq.det.eq.2}) yields%
\begin{align*}
\det A  &  =\underbrace{\sum_{\sigma\in S_{n}}}_{\substack{=\sum
_{P\subseteq\left\{  1,2,\ldots,n\right\}  }\sum_{\substack{\sigma\in
S_{n};\\\sigma^{-1}\left(  Q\right)  =P}}\\\text{(because for every }\sigma\in
S_{n}\text{, the set}\\\sigma^{-1}\left(  Q\right)  \text{ is a subset of
}\left\{  1,2,\ldots,n\right\}  \text{)}}}\left(  -1\right)  ^{\sigma}%
\prod_{i=1}^{n}a_{i,\sigma\left(  i\right)  }\\
&  =\sum_{P\subseteq\left\{  1,2,\ldots,n\right\}  }\underbrace{\sum
_{\substack{\sigma\in S_{n};\\\sigma^{-1}\left(  Q\right)  =P}}}%
_{\substack{=\sum_{\substack{\sigma\in S_{n};\\\sigma\left(  P\right)
=Q}}\\\text{(because for every }\sigma\in S_{n}\text{,}\\\text{the statement
}\left(  \sigma^{-1}\left(  Q\right)  =P\right)  \\\text{is equivalent
to}\\\text{the statement }\left(  \sigma\left(  P\right)  =Q\right)
\\\text{(by Lemma \ref{lem.sol.det.laplace-multi.4} \textbf{(b)}))}}}\left(
-1\right)  ^{\sigma}\prod_{i=1}^{n}a_{i,\sigma\left(  i\right)  }\\
&  =\sum_{P\subseteq\left\{  1,2,\ldots,n\right\}  }\sum_{\substack{\sigma\in
S_{n};\\\sigma\left(  P\right)  =Q}}\left(  -1\right)  ^{\sigma}\prod
_{i=1}^{n}a_{i,\sigma\left(  i\right)  }\\
&  =\sum_{\substack{P\subseteq\left\{  1,2,\ldots,n\right\}  ;\\\left\vert
P\right\vert =\left\vert Q\right\vert }}\sum_{\substack{\sigma\in
S_{n};\\\sigma\left(  P\right)  =Q}}\left(  -1\right)  ^{\sigma}\prod
_{i=1}^{n}a_{i,\sigma\left(  i\right)  }+\sum_{\substack{P\subseteq\left\{
1,2,\ldots,n\right\}  ;\\\left\vert P\right\vert \neq\left\vert Q\right\vert
}}\underbrace{\sum_{\substack{\sigma\in S_{n};\\\sigma\left(  P\right)
=Q}}\left(  -1\right)  ^{\sigma}\prod_{i=1}^{n}a_{i,\sigma\left(  i\right)  }%
}_{\substack{=0\\\text{(by (\ref{pf.thm.det.laplace-multi.emptysum}))}}}\\
&  \ \ \ \ \ \ \ \ \ \ \left(
\begin{array}
[c]{c}%
\text{since every subset }P\text{ of }\left\{  1,2,\ldots,n\right\}  \text{
satisfies}\\
\text{either }\left\vert P\right\vert =\left\vert Q\right\vert \text{ or
}\left\vert P\right\vert \neq\left\vert Q\right\vert \text{ (but not both)}%
\end{array}
\right) \\
&  =\sum_{\substack{P\subseteq\left\{  1,2,\ldots,n\right\}  ;\\\left\vert
P\right\vert =\left\vert Q\right\vert }}\sum_{\substack{\sigma\in
S_{n};\\\sigma\left(  P\right)  =Q}}\left(  -1\right)  ^{\sigma}\prod
_{i=1}^{n}a_{i,\sigma\left(  i\right)  }+\underbrace{\sum
_{\substack{P\subseteq\left\{  1,2,\ldots,n\right\}  ;\\\left\vert
P\right\vert \neq\left\vert Q\right\vert }}0}_{=0}\\
&  =\sum_{\substack{P\subseteq\left\{  1,2,\ldots,n\right\}  ;\\\left\vert
P\right\vert =\left\vert Q\right\vert }}\underbrace{\sum_{\substack{\sigma\in
S_{n};\\\sigma\left(  P\right)  =Q}}\left(  -1\right)  ^{\sigma}\prod
_{i=1}^{n}a_{i,\sigma\left(  i\right)  }}_{\substack{=\left(  -1\right)
^{\sum P+\sum Q}\det\left(  \operatorname*{sub}\nolimits_{w\left(  P\right)
}^{w\left(  Q\right)  }A\right)  \det\left(  \operatorname*{sub}%
\nolimits_{w\left(  \widetilde{P}\right)  }^{w\left(  \widetilde{Q}\right)
}A\right)  \\\text{(by Lemma \ref{lem.sol.det.laplace-multi.2})}}}\\
&  =\sum_{\substack{P\subseteq\left\{  1,2,\ldots,n\right\}  ;\\\left\vert
P\right\vert =\left\vert Q\right\vert }}\left(  -1\right)  ^{\sum P+\sum
Q}\det\left(  \operatorname*{sub}\nolimits_{w\left(  P\right)  }^{w\left(
Q\right)  }A\right)  \det\left(  \operatorname*{sub}\nolimits_{w\left(
\widetilde{P}\right)  }^{w\left(  \widetilde{Q}\right)  }A\right)  .
\end{align*}
This proves Theorem \ref{thm.det.laplace-multi} \textbf{(b)}.
\end{proof}

\begin{proof}
[Solution to Exercise \ref{exe.det.laplace-multi}.]We have now proven both
Lemma \ref{lem.det.laplace-multi.Apq} and Theorem \ref{thm.det.laplace-multi}.
Thus, Exercise \ref{exe.det.laplace-multi} is solved.
\end{proof}

\subsection{Solution to Exercise \ref{exe.det.laplace-multi.0}}

Throughout this section, we shall use the notations introduced in Definition
\ref{def.submatrix} and in Definition \ref{def.sect.laplace.notations}.

Before we start solving Exercise \ref{exe.det.laplace-multi.0}, let us show a
trivial lemma:

\begin{lemma}
\label{lem.sol.det.laplace-multi.0.subAT}Let $n\in\mathbb{N}$ and
$m\in\mathbb{N}$. Let $A$ be an $n\times m$-matrix. Let $\mathbf{i}$ be a
finite list of elements of $\left\{  1,2,\ldots,n\right\}  $. Let $\mathbf{j}$
be a finite list of elements of $\left\{  1,2,\ldots,m\right\}  $. Then,
$\left(  \operatorname*{sub}\nolimits_{\mathbf{i}}^{\mathbf{j}}A\right)
^{T}=\operatorname*{sub}\nolimits_{\mathbf{j}}^{\mathbf{i}}\left(
A^{T}\right)  $.
\end{lemma}

\begin{vershort}
\begin{proof}
[Proof of Lemma \ref{lem.sol.det.laplace-multi.0.subAT}.]Write the list
$\mathbf{i}$ as $\left(  i_{1},i_{2},\ldots,i_{u}\right)  $. Write the list
$\mathbf{j}$ as $\left(  j_{1},j_{2},\ldots,j_{v}\right)  $. Then, Lemma
\ref{lem.sol.det.laplace-multi.0.subAT} follows immediately from Proposition
\ref{prop.submatrix.easy} \textbf{(e)}.
\end{proof}
\end{vershort}

\begin{verlong}
\begin{proof}
[Proof of Lemma \ref{lem.sol.det.laplace-multi.0.subAT}.]Write the list
$\mathbf{i}$ in the form $\mathbf{i}=\left(  i_{1},i_{2},\ldots,i_{u}\right)
$. Then, $\left(  i_{1},i_{2},\ldots,i_{u}\right)  =\mathbf{i}$ is a list of
elements of $\left\{  1,2,\ldots,n\right\}  $. In other words, $i_{1}%
,i_{2},\ldots,i_{u}$ are some elements of $\left\{  1,2,\ldots,n\right\}  $.

Write the list $\mathbf{j}$ in the form $\mathbf{j}=\left(  j_{1},j_{2}%
,\ldots,j_{v}\right)  $. Then, $\left(  j_{1},j_{2},\ldots,j_{v}\right)
=\mathbf{j}$ is a list of elements of $\left\{  1,2,\ldots,m\right\}  $. In
other words, $j_{1},j_{2},\ldots,j_{v}$ are some elements of $\left\{
1,2,\ldots,m\right\}  $.

From $\mathbf{i}=\left(  i_{1},i_{2},\ldots,i_{u}\right)  $ and $\mathbf{j}%
=\left(  j_{1},j_{2},\ldots,j_{v}\right)  $, we obtain%
\[
\operatorname*{sub}\nolimits_{\mathbf{i}}^{\mathbf{j}}A=\operatorname*{sub}%
\nolimits_{\left(  i_{1},i_{2},\ldots,i_{u}\right)  }^{\left(  j_{1}%
,j_{2},\ldots,j_{v}\right)  }A=\operatorname*{sub}\nolimits_{i_{1}%
,i_{2},\ldots,i_{u}}^{j_{1},j_{2},\ldots,j_{v}}A.
\]
Hence,%
\begin{equation}
\left(  \underbrace{\operatorname*{sub}\nolimits_{\mathbf{i}}^{\mathbf{j}}%
A}_{=\operatorname*{sub}\nolimits_{i_{1},i_{2},\ldots,i_{u}}^{j_{1}%
,j_{2},\ldots,j_{v}}A}\right)  ^{T}=\left(  \operatorname*{sub}%
\nolimits_{i_{1},i_{2},\ldots,i_{u}}^{j_{1},j_{2},\ldots,j_{v}}A\right)
^{T}=\operatorname*{sub}\nolimits_{j_{1},j_{2},\ldots,j_{v}}^{i_{1}%
,i_{2},\ldots,i_{u}}\left(  A^{T}\right)
\label{pf.lem.sol.det.laplace-multi.0.subAT.1}%
\end{equation}
(by Proposition \ref{prop.submatrix.easy} \textbf{(e)}). From $\mathbf{j}%
=\left(  j_{1},j_{2},\ldots,j_{v}\right)  $ and $\mathbf{i}=\left(
i_{1},i_{2},\ldots,i_{u}\right)  $, we obtain%
\[
\operatorname*{sub}\nolimits_{\mathbf{j}}^{\mathbf{i}}\left(  A^{T}\right)
=\operatorname*{sub}\nolimits_{\left(  j_{1},j_{2},\ldots,j_{v}\right)
}^{\left(  i_{1},i_{2},\ldots,i_{u}\right)  }\left(  A^{T}\right)
=\operatorname*{sub}\nolimits_{j_{1},j_{2},\ldots,j_{v}}^{i_{1},i_{2}%
,\ldots,i_{u}}\left(  A^{T}\right)  .
\]
Comparing this with (\ref{pf.lem.sol.det.laplace-multi.0.subAT.1}), we obtain
$\left(  \operatorname*{sub}\nolimits_{\mathbf{i}}^{\mathbf{j}}A\right)
^{T}=\operatorname*{sub}\nolimits_{\mathbf{j}}^{\mathbf{i}}\left(
A^{T}\right)  $. This proves Lemma \ref{lem.sol.det.laplace-multi.0.subAT}.
\end{proof}
\end{verlong}

\begin{corollary}
\label{cor.sol.det.laplace-multi.0.detsubAT}Let $n\in\mathbb{N}$ and
$m\in\mathbb{N}$.

Let $A$ be an $n\times m$-matrix. Let $U$ be a subset of $\left\{
1,2,\ldots,m\right\}  $. Let $V$ be a subset of $\left\{  1,2,\ldots
,n\right\}  $. Assume that $\left\vert U\right\vert =\left\vert V\right\vert
$. Then,%
\[
\det\left(  \operatorname*{sub}\nolimits_{w\left(  U\right)  }^{w\left(
V\right)  }\left(  A^{T}\right)  \right)  =\det\left(  \operatorname*{sub}%
\nolimits_{w\left(  V\right)  }^{w\left(  U\right)  }A\right)  .
\]

\end{corollary}

\begin{vershort}
\begin{proof}
[Proof of Corollary \ref{cor.sol.det.laplace-multi.0.detsubAT}.]Let
$k=\left\vert U\right\vert =\left\vert V\right\vert $. Then, each of the lists
$w\left(  U\right)  $ and $w\left(  V\right)  $ is a list of $k$ elements.
Hence, $\operatorname*{sub}\nolimits_{w\left(  V\right)  }^{w\left(  U\right)
}A$ is a $k\times k$-matrix. Thus, Exercise \ref{exe.ps4.4} (applied to $k$
and $\operatorname*{sub}\nolimits_{w\left(  V\right)  }^{w\left(  U\right)
}A$ instead of $n$ and $A$) yields%
\begin{equation}
\det\left(  \left(  \operatorname*{sub}\nolimits_{w\left(  V\right)
}^{w\left(  U\right)  }A\right)  ^{T}\right)  =\det\left(  \operatorname*{sub}%
\nolimits_{w\left(  V\right)  }^{w\left(  U\right)  }A\right)  .
\label{pf.cor.sol.det.laplace-multi.0.detsubAT.short.1}%
\end{equation}
But $w\left(  U\right)  $ is a list of elements of $\left\{  1,2,\ldots
,m\right\}  $ (since $U$ is a subset of $\left\{  1,2,\ldots,m\right\}  $).
Similarly, $w\left(  V\right)  $ is a list of elements of $\left\{
1,2,\ldots,n\right\}  $. Hence, Lemma \ref{lem.sol.det.laplace-multi.0.subAT}
(applied to $\mathbf{i}=w\left(  V\right)  $ and $\mathbf{j}=w\left(
U\right)  $) yields $\left(  \operatorname*{sub}\nolimits_{w\left(  V\right)
}^{w\left(  U\right)  }A\right)  ^{T}=\operatorname*{sub}\nolimits_{w\left(
U\right)  }^{w\left(  V\right)  }\left(  A^{T}\right)  $. Taking determinants
on both sides of this equality, we obtain%
\[
\det\left(  \left(  \operatorname*{sub}\nolimits_{w\left(  V\right)
}^{w\left(  U\right)  }A\right)  ^{T}\right)  =\det\left(  \operatorname*{sub}%
\nolimits_{w\left(  U\right)  }^{w\left(  V\right)  }\left(  A^{T}\right)
\right)  .
\]
Comparing this with (\ref{pf.cor.sol.det.laplace-multi.0.detsubAT.short.1}),
we obtain $\det\left(  \operatorname*{sub}\nolimits_{w\left(  U\right)
}^{w\left(  V\right)  }\left(  A^{T}\right)  \right)  =\det\left(
\operatorname*{sub}\nolimits_{w\left(  V\right)  }^{w\left(  U\right)
}A\right)  $. This proves Corollary \ref{cor.sol.det.laplace-multi.0.detsubAT}.
\end{proof}
\end{vershort}

\begin{verlong}
\begin{proof}
[Proof of Corollary \ref{cor.sol.det.laplace-multi.0.detsubAT}.]We have
$U\subseteq\left\{  1,2,\ldots,m\right\}  $. Thus, the set $U$ is finite
(since the set $\left\{  1,2,\ldots,m\right\}  $ is finite). Hence,
$\left\vert U\right\vert \in\mathbb{N}$. Define $k\in\mathbb{N}$ by
$k=\left\vert U\right\vert $. Thus, $k=\left\vert U\right\vert =\left\vert
V\right\vert $.

Recall that $w\left(  U\right)  $ is the list of all elements of $U$ in
increasing order (with no repetitions) (by the definition of $w\left(
U\right)  $). Thus, $w\left(  U\right)  $ is a list of size $\left\vert
U\right\vert $. In other words, $w\left(  U\right)  $ is a list of size $k$
(since $\left\vert U\right\vert =k$).

Write the list $w\left(  U\right)  $ in the form $w\left(  U\right)  =\left(
i_{1},i_{2},\ldots,i_{k}\right)  $. (This is possible, since $w\left(
U\right)  $ is a list of size $k$.)

Recall that $w\left(  V\right)  $ is the list of all elements of $V$ in
increasing order (with no repetitions) (by the definition of $w\left(
V\right)  $). Thus, $w\left(  V\right)  $ is a list of size $\left\vert
V\right\vert $. In other words, $w\left(  V\right)  $ is a list of size $k$
(since $\left\vert V\right\vert =k$).

Write the list $w\left(  V\right)  $ in the form $w\left(  V\right)  =\left(
j_{1},j_{2},\ldots,j_{k}\right)  $. (This is possible, since $w\left(
V\right)  $ is a list of size $k$.)

From $w\left(  U\right)  =\left(  i_{1},i_{2},\ldots,i_{k}\right)  $ and
$w\left(  V\right)  =\left(  j_{1},j_{2},\ldots,j_{k}\right)  $, we obtain%
\[
\operatorname*{sub}\nolimits_{w\left(  V\right)  }^{w\left(  U\right)
}A=\operatorname*{sub}\nolimits_{\left(  j_{1},j_{2},\ldots,j_{k}\right)
}^{\left(  i_{1},i_{2},\ldots,i_{k}\right)  }A=\operatorname*{sub}%
\nolimits_{j_{1},j_{2},\ldots,j_{k}}^{i_{1},i_{2},\ldots,i_{k}}A.
\]
Hence, $\operatorname*{sub}\nolimits_{w\left(  V\right)  }^{w\left(  U\right)
}A$ is a $k\times k$-matrix (since $\operatorname*{sub}\nolimits_{j_{1}%
,j_{2},\ldots,j_{k}}^{i_{1},i_{2},\ldots,i_{k}}A$ is a $k\times k$-matrix).
Thus, Exercise \ref{exe.ps4.4} (applied to $k$ and $\operatorname*{sub}%
\nolimits_{w\left(  V\right)  }^{w\left(  U\right)  }A$ instead of $n$ and
$A$) yields%
\begin{equation}
\det\left(  \left(  \operatorname*{sub}\nolimits_{w\left(  V\right)
}^{w\left(  U\right)  }A\right)  ^{T}\right)  =\det\left(  \operatorname*{sub}%
\nolimits_{w\left(  V\right)  }^{w\left(  U\right)  }A\right)  .
\label{pf.cor.sol.det.laplace-multi.0.detsubAT.det1}%
\end{equation}

We know that $w\left(  U\right)  $ is the list of all elements of $U$ in
increasing order (with no repetitions). Thus, $w\left(  U\right)  $ is a list
of elements of $U$. In other words, the entries of $w\left(  U\right)  $ are
elements of $U$. Hence, the entries of $w\left(  U\right)  $ are elements of
$\left\{  1,2,\ldots,m\right\}  $ (since every element of $U$ is an element of
$\left\{  1,2,\ldots,m\right\}  $ (since $U\subseteq\left\{  1,2,\ldots
,m\right\}  $)). In other words, $w\left(  U\right)  $ is a list of elements
of $\left\{  1,2,\ldots,m\right\}  $. The same argument (applied to $V$ and
$n$ instead of $U$ and $m$) shows that $w\left(  V\right)  $ is a list of
elements of $\left\{  1,2,\ldots,n\right\}  $. Hence, Lemma
\ref{lem.sol.det.laplace-multi.0.subAT} (applied to $\mathbf{i}=w\left(
V\right)  $ and $\mathbf{j}=w\left(  U\right)  $) yields $\left(
\operatorname*{sub}\nolimits_{w\left(  V\right)  }^{w\left(  U\right)
}A\right)  ^{T}=\operatorname*{sub}\nolimits_{w\left(  U\right)  }^{w\left(
V\right)  }\left(  A^{T}\right)  $. Therefore,%
\[
\det\left(  \underbrace{\left(  \operatorname*{sub}\nolimits_{w\left(
V\right)  }^{w\left(  U\right)  }A\right)  ^{T}}_{=\operatorname*{sub}%
\nolimits_{w\left(  U\right)  }^{w\left(  V\right)  }\left(  A^{T}\right)
}\right)  =\det\left(  \operatorname*{sub}\nolimits_{w\left(  U\right)
}^{w\left(  V\right)  }\left(  A^{T}\right)  \right)  .
\]
Comparing this with (\ref{pf.cor.sol.det.laplace-multi.0.detsubAT.det1}), we
obtain $\det\left(  \operatorname*{sub}\nolimits_{w\left(  U\right)
}^{w\left(  V\right)  }\left(  A^{T}\right)  \right)  =\det\left(
\operatorname*{sub}\nolimits_{w\left(  V\right)  }^{w\left(  U\right)
}A\right)  $. This proves Corollary \ref{cor.sol.det.laplace-multi.0.detsubAT}.
\end{proof}
\end{verlong}

\begin{proof}
[Solution to Exercise \ref{exe.det.laplace-multi.0}.]\textbf{(a)} Let $P$ be a
subset of $\left\{  1,2,\ldots,n\right\}  $ satisfying $\left\vert
P\right\vert =\left\vert R\right\vert $ and $P\neq R$.

\begin{vershort}
Let $k=\left\vert P\right\vert $. The definition of $\widetilde{P}$ yields
$\widetilde{P}=\left\{  1,2,\ldots,n\right\}  \setminus P$. Since $P$ is a
subset of $\left\{  1,2,\ldots,n\right\}  $, this leads to $\left\vert
\widetilde{P}\right\vert =\underbrace{\left\vert \left\{  1,2,\ldots
,n\right\}  \right\vert }_{=n}-\underbrace{\left\vert P\right\vert }_{=k}%
=n-k$. Thus, $n-k=\left\vert \widetilde{P}\right\vert \geq0$, so that $n\geq
k$ and thus $k\in\left\{  0,1,\ldots,n\right\}  $.
\end{vershort}

\begin{verlong}
The definition of $\widetilde{P}$ yields $\widetilde{P}=\left\{
1,2,\ldots,n\right\}  \setminus P\subseteq\left\{  1,2,\ldots,n\right\}  $.
Thus, $\widetilde{P}$ is a finite set (since $\left\{  1,2,\ldots,n\right\}  $
is a finite set). Hence, $\left\vert \widetilde{P}\right\vert \in\mathbb{N}$.

Also, $P\subseteq\left\{  1,2,\ldots,n\right\}  $. Thus, $P$ is a finite set
(since $\left\{  1,2,\ldots,n\right\}  $ is a finite set). Hence, $\left\vert
P\right\vert \in\mathbb{N}$.

Define $k\in\mathbb{N}$ by $k=\left\vert P\right\vert $. (This is
well-defined, since $\left\vert P\right\vert \in\mathbb{N}$.)

Also,%
\begin{align*}
\left\vert \underbrace{\widetilde{P}}_{=\left\{  1,2,\ldots,n\right\}
\setminus P}\right\vert  &  =\left\vert \left\{  1,2,\ldots,n\right\}
\setminus P\right\vert =\underbrace{\left\vert \left\{  1,2,\ldots,n\right\}
\right\vert }_{=n}-\underbrace{\left\vert P\right\vert }_{=k}%
\ \ \ \ \ \ \ \ \ \ \left(  \text{since }P\subseteq\left\{  1,2,\ldots
,n\right\}  \right) \\
&  =n-k.
\end{align*}
Hence, $n-k=\left\vert \widetilde{P}\right\vert \in\mathbb{N}$; thus,
$n-k\geq0$. In other words, $k\leq n$. Combined with $k\geq0$ (since
$k=\left\vert P\right\vert \in\mathbb{N}$), this yields $k\in\left\{
0,1,\ldots,n\right\}  $.

We have $k\in\mathbb{N}$, so that $k\geq0$. Thus, $\underbrace{k}_{\geq
0}+1\geq1$.
\end{verlong}

Let $\left[  n\right]  $ denote the set $\left\{  1,2,\ldots,n\right\}  $.
Thus, $\left[  n\right]  =\left\{  1,2,\ldots,n\right\}  $.

We have $k=\left\vert P\right\vert =\left\vert R\right\vert $, so that
$\left\vert R\right\vert =k$.

From $P\neq R$, we can easily deduce that $\widetilde{P}\cap R\neq\varnothing
$\ \ \ \ \footnote{\textit{Proof.} Assume the contrary (for the sake of
contradiction). Thus, $\widetilde{P}\cap R=\varnothing$.
\par
Any three sets $X$, $Y$ and $Z$ satisfy $\left(  X\cap Y\right)  \setminus
Z=X\cap\left(  Y\setminus Z\right)  $. Applying this to $X=R$, $Y=\left\{
1,2,\ldots,n\right\}  $ and $Z=P$, we obtain%
\[
\left(  R\cap\left\{  1,2,\ldots,n\right\}  \right)  \setminus P=R\cap
\underbrace{\left(  \left\{  1,2,\ldots,n\right\}  \setminus P\right)
}_{=\widetilde{P}}=R\cap\widetilde{P}=\widetilde{P}\cap R=\varnothing.
\]
Comparing this with $\underbrace{\left(  R\cap\left\{  1,2,\ldots,n\right\}
\right)  }_{\substack{=R\\\text{(since }R\subseteq\left\{  1,2,\ldots
,n\right\}  \text{)}}}\setminus P=R\setminus P$, we obtain $R\setminus
P=\varnothing$. In other words, $R\subseteq P$. Combining this with $R\neq P$
(since $P\neq R$), we conclude that $R$ is a proper subset of $P$.
\par
But $P$ is a finite set. Hence, every proper subset of $P$ has a size strictly
smaller than $\left\vert P\right\vert $. In other words, if $Y$ is a proper
subset of $P$, then $\left\vert Y\right\vert <\left\vert P\right\vert $.
Applying this to $Y=R$, we conclude that $\left\vert R\right\vert <\left\vert
P\right\vert $. Hence, $\left\vert R\right\vert <\left\vert P\right\vert
=\left\vert R\right\vert $. This is absurd. This contradiction proves that our
assumption was wrong, qed.}.

Let $\left(  r_{1},r_{2},\ldots,r_{k}\right)  $ be the list of all elements of
$R$ in increasing order (with no repetitions). (This is well-defined, because
$\left\vert R\right\vert =k$.)

Let $\left(  t_{1},t_{2},\ldots,t_{n-k}\right)  $ be the list of all elements
of $\widetilde{P}$ in increasing order (with no repetitions). (This is
well-defined, because $\left\vert \widetilde{P}\right\vert =n-k$.)

\begin{vershort}
We know that $\left(  r_{1},r_{2},\ldots,r_{k}\right)  $ is a list of elements
of $R$. Hence, the elements $r_{1},r_{2},\ldots,r_{k}$ belong to $R$, and thus
to $\left[  n\right]  $ (since $R\subseteq\left\{  1,2,\ldots,n\right\}
=\left[  n\right]  $).

We know that $\left(  t_{1},t_{2},\ldots,t_{n-k}\right)  $ is a list of
elements of $\widetilde{P}$. Hence, the elements $t_{1},t_{2},\ldots,t_{n-k}$
belong to $\widetilde{P}$, and thus to $\left[  n\right]  $ (since
$\widetilde{P}\subseteq\left\{  1,2,\ldots,n\right\}  =\left[  n\right]  $).
\end{vershort}

\begin{verlong}
We know that $\left(  r_{1},r_{2},\ldots,r_{k}\right)  $ is a list of elements
of $R$ (since $\left(  r_{1},r_{2},\ldots,r_{k}\right)  $ is the list of all
elements of $R$ in increasing order (with no repetitions)). Hence, the $k$
elements $r_{1},r_{2},\ldots,r_{k}$ belong to $R$, and thus belong to $\left[
n\right]  $ as well (since $R\subseteq\left\{  1,2,\ldots,n\right\}  =\left[
n\right]  $).

We know that $\left(  t_{1},t_{2},\ldots,t_{n-k}\right)  $ is a list of
elements of $\widetilde{P}$ (since $\left(  t_{1},t_{2},\ldots,t_{n-k}\right)
$ is the list of all elements of $\widetilde{P}$ in increasing order (with no
repetitions)). Hence, the $n-k$ elements $t_{1},t_{2},\ldots,t_{n-k}$ belong
to $\widetilde{P}$, and thus belong to $\left[  n\right]  $ as well (since
$\widetilde{P}\subseteq\left\{  1,2,\ldots,n\right\}  =\left[  n\right]  $).
\end{verlong}

\begin{vershort}
The $n$ elements $r_{1},r_{2},\ldots,r_{k},t_{1},t_{2},\ldots,t_{n-k}$ belong
to $\left[  n\right]  $ (since the $k$ elements $r_{1},r_{2},\ldots,r_{k}$
belong to $\left[  n\right]  $, and since the $n-k$ elements $t_{1}%
,t_{2},\ldots,t_{n-k}$ belong to $\left[  n\right]  $). Hence, we can define
an $n$-tuple $\left(  \kappa_{1},\kappa_{2},\ldots,\kappa_{n}\right)
\in\left[  n\right]  ^{n}$ by%
\[
\left(  \kappa_{1},\kappa_{2},\ldots,\kappa_{n}\right)  =\left(  r_{1}%
,r_{2},\ldots,r_{k},t_{1},t_{2},\ldots,t_{n-k}\right)  .
\]
Consider this $n$-tuple. Due to its definition, we have%
\begin{equation}
\left(  \kappa_{i}=r_{i}\ \ \ \ \ \ \ \ \ \ \text{for every }i\in\left\{
1,2,\ldots,k\right\}  \right)  \label{sol.det.laplace-multi.0.short.a.kappai1}%
\end{equation}
and%
\begin{equation}
\left(  \kappa_{i}=t_{i-k}\ \ \ \ \ \ \ \ \ \ \text{for every }i\in\left\{
k+1,k+2,\ldots,n\right\}  \right)  .
\label{sol.det.laplace-multi.0.short.a.kappai2}%
\end{equation}

Define a map $\kappa:\left[  n\right]  \rightarrow\left[  n\right]  $ by%
\[
\left(  \kappa\left(  i\right)  =\kappa_{i}\text{ for every }i\in\left[
n\right]  \right)  .
\]
(This makes sense, since $\kappa_{i}\in\left[  n\right]  $ for every
$i\in\left[  n\right]  $.) Then, every $i\in\left\{  1,2,\ldots,k\right\}  $
satisfies%
\begin{align}
\kappa\left(  i\right)   &  =\kappa_{i}\ \ \ \ \ \ \ \ \ \ \left(  \text{by
the definition of }\kappa\right) \nonumber\\
&  =r_{i}\ \ \ \ \ \ \ \ \ \ \left(  \text{by
(\ref{sol.det.laplace-multi.0.short.a.kappai1})}\right)  .
\label{sol.det.laplace-multi.0.short.a.kappa1}%
\end{align}
Also, every $j\in\left\{  1,2,\ldots,n-k\right\}  $ satisfies%
\begin{align}
\kappa\left(  k+j\right)   &  =\kappa_{k+j}\ \ \ \ \ \ \ \ \ \ \left(
\text{by the definition of }\kappa\right) \nonumber\\
&  =t_{\left(  k+j\right)  -k}\ \ \ \ \ \ \ \ \ \ \left(  \text{by
(\ref{sol.det.laplace-multi.0.short.a.kappai2}), applied to }i=k+j\right)
\nonumber\\
&  =t_{j}. \label{sol.det.laplace-multi.0.short.a.kappa2}%
\end{align}

\end{vershort}

\begin{verlong}
The $k+\left(  n-k\right)  $ elements $r_{1},r_{2},\ldots,r_{k},t_{1}%
,t_{2},\ldots,t_{n-k}$ belong to $\left[  n\right]  $ (since the $k$ elements
$r_{1},r_{2},\ldots,r_{k}$ belong to $\left[  n\right]  $, and since the $n-k$
elements $t_{1},t_{2},\ldots,t_{n-k}$ belong to $\left[  n\right]  $). In
other words, $\left(  r_{1},r_{2},\ldots,r_{k},t_{1},t_{2},\ldots
,t_{n-k}\right)  \in\left[  n\right]  ^{k+\left(  n-k\right)  }$. Thus,%
\[
\left(  r_{1},r_{2},\ldots,r_{k},t_{1},t_{2},\ldots,t_{n-k}\right)  \in\left[
n\right]  ^{k+\left(  n-k\right)  }=\left[  n\right]  ^{n}%
\]
(since $k+\left(  n-k\right)  =n$). In other words, $\left(  r_{1}%
,r_{2},\ldots,r_{k},t_{1},t_{2},\ldots,t_{n-k}\right)  $ is an $n$-tuple of
elements of $\left[  n\right]  $. Write this $n$-tuple $\left(  r_{1}%
,r_{2},\ldots,r_{k},t_{1},t_{2},\ldots,t_{n-k}\right)  $ in the form $\left(
\kappa_{1},\kappa_{2},\ldots,\kappa_{n}\right)  $. Thus, $\left(  \kappa
_{1},\kappa_{2},\ldots,\kappa_{n}\right)  =\left(  r_{1},r_{2},\ldots
,r_{k},t_{1},t_{2},\ldots,t_{n-k}\right)  $. In other words,%
\begin{equation}
\kappa_{i}=%
\begin{cases}
r_{i}, & \text{if }i\leq k;\\
t_{i-k}, & \text{if }i>k
\end{cases}
\ \ \ \ \ \ \ \ \ \ \text{for every }i\in\left\{  1,2,\ldots,n\right\}  .
\label{sol.det.laplace-multi.0.a.def-kappa}%
\end{equation}
Moreover, $\left(  \kappa_{1},\kappa_{2},\ldots,\kappa_{n}\right)  \in\left[
n\right]  ^{n}$ (since $\left(  \kappa_{1},\kappa_{2},\ldots,\kappa
_{n}\right)  $ is an $n$-tuple of elements of $\left[  n\right]  $). Hence,
$\kappa_{i}\in\left[  n\right]  $ for each $i\in\left\{  1,2,\ldots,n\right\}
$. In other words, $\kappa_{i}\in\left[  n\right]  $ for each $i\in\left[
n\right]  $ (since $\left[  n\right]  =\left\{  1,2,\ldots,n\right\}  $).
Thus, we can define a map $\kappa:\left[  n\right]  \rightarrow\left[
n\right]  $ by%
\[
\left(  \kappa\left(  i\right)  =\kappa_{i}\text{ for every }i\in\left[
n\right]  \right)  .
\]
Consider this $\kappa$.

We have%
\begin{equation}
\kappa\left(  i\right)  =r_{i}\ \ \ \ \ \ \ \ \ \ \text{for every }%
i\in\left\{  1,2,\ldots,k\right\}  . \label{sol.det.laplace-multi.0.a.kappa1}%
\end{equation}

[\textit{Proof of (\ref{sol.det.laplace-multi.0.a.kappa1}):} Let $i\in\left\{
1,2,\ldots,k\right\}  $. Then, $i\in\left\{  1,2,\ldots,k\right\}
\subseteq\left\{  1,2,\ldots,n\right\}  $ (since $k\leq n$), so that
$i\in\left\{  1,2,\ldots,n\right\}  =\left[  n\right]  $. Hence,
$\kappa\left(  i\right)  =\kappa_{i}$ (by the definition of $\kappa$). Thus,%
\begin{align*}
\kappa\left(  i\right)   &  =\kappa_{i}=%
\begin{cases}
r_{i}, & \text{if }i\leq k;\\
t_{i-k}, & \text{if }i>k
\end{cases}
\ \ \ \ \ \ \ \ \ \ \left(  \text{by the definition of }\kappa_{i}\right) \\
&  =r_{i}\ \ \ \ \ \ \ \ \ \ \left(  \text{since }i\leq k\text{ (since }%
i\in\left\{  1,2,\ldots,k\right\}  \text{)}\right)  .
\end{align*}
This proves (\ref{sol.det.laplace-multi.0.a.kappa1}).]

Furthermore,%
\begin{equation}
\kappa\left(  k+i\right)  =t_{i}\ \ \ \ \ \ \ \ \ \ \text{for every }%
i\in\left\{  1,2,\ldots,n-k\right\}  .
\label{sol.det.laplace-multi.0.a.kappa2}%
\end{equation}

[\textit{Proof of (\ref{sol.det.laplace-multi.0.a.kappa2}):} Let $i\in\left\{
1,2,\ldots,n-k\right\}  $. Thus, $i\geq1>0$. But $i\in\left\{  1,2,\ldots
,n-k\right\}  $, and thus%
\begin{align*}
k+i  &  \in\left\{  k+1,k+2,\ldots,k+\left(  n-k\right)  \right\}  =\left\{
k+1,k+2,\ldots,n\right\}  \ \ \ \ \ \ \ \ \ \ \left(  \text{since }k+\left(
n-k\right)  =n\right) \\
&  \subseteq\left\{  1,2,\ldots,n\right\}  \ \ \ \ \ \ \ \ \ \ \left(
\text{since }k+1\geq1\right) \\
&  =\left[  n\right]  .
\end{align*}
Hence, $\kappa\left(  k+i\right)  =\kappa_{k+i}$ (by the definition of
$\kappa$). Hence,%
\begin{align*}
\kappa\left(  k+i\right)   &  =\kappa_{k+i}=%
\begin{cases}
r_{k+i}, & \text{if }k+i\leq k;\\
t_{\left(  k+i\right)  -k}, & \text{if }k+i>k
\end{cases}
\ \ \ \ \ \ \ \ \ \ \left(  \text{by the definition of }\kappa_{k+i}\right) \\
&  =t_{\left(  k+i\right)  -k}\ \ \ \ \ \ \ \ \ \ \left(  \text{since
}k+\underbrace{i}_{>0}>k\right) \\
&  =t_{i}\ \ \ \ \ \ \ \ \ \ \left(  \text{since }\left(  k+i\right)
-k=i\right)  .
\end{align*}
This proves (\ref{sol.det.laplace-multi.0.a.kappa2}).]
\end{verlong}

\begin{vershort}
We have%
\begin{equation}
w\left(  R\right)  =\left(  \kappa\left(  1\right)  ,\kappa\left(  2\right)
,\ldots,\kappa\left(  k\right)  \right)
\label{sol.det.laplace-multi.0.short.a.wR=}%
\end{equation}
\footnote{\textit{Proof of (\ref{sol.det.laplace-multi.0.short.a.wR=}):} The
lists $w\left(  R\right)  $ and $\left(  r_{1},r_{2},\ldots,r_{k}\right)  $
must be identical (since each of them is the list of all elements of $R$ in
increasing order (with no repetitions)). In other words, $w\left(  R\right)
=\left(  r_{1},r_{2},\ldots,r_{k}\right)  $. Now,
(\ref{sol.det.laplace-multi.0.short.a.kappa1}) shows that $\kappa\left(
i\right)  =r_{i}$ for every $i\in\left\{  1,2,\ldots,k\right\}  $. In other
words, $\left(  \kappa\left(  1\right)  ,\kappa\left(  2\right)
,\ldots,\kappa\left(  k\right)  \right)  =\left(  r_{1},r_{2},\ldots
,r_{k}\right)  $. Comparing this with $w\left(  R\right)  =\left(  r_{1}%
,r_{2},\ldots,r_{k}\right)  $, we obtain $w\left(  R\right)  =\left(
\kappa\left(  1\right)  ,\kappa\left(  2\right)  ,\ldots,\kappa\left(
k\right)  \right)  $. This proves (\ref{sol.det.laplace-multi.0.short.a.wR=}%
).} and%
\begin{equation}
w\left(  \widetilde{P}\right)  =\left(  \kappa\left(  k+1\right)
,\kappa\left(  k+2\right)  ,\ldots,\kappa\left(  n\right)  \right)
\label{sol.det.laplace-multi.0.short.a.wPs=}%
\end{equation}
\footnote{\textit{Proof of (\ref{sol.det.laplace-multi.0.short.a.wPs=}):} The
lists $w\left(  \widetilde{P}\right)  $ and $\left(  t_{1},t_{2}%
,\ldots,t_{n-k}\right)  $ must be identical (since each of them is the list of
all elements of $\widetilde{P}$ in increasing order (with no repetitions)). In
other words, $w\left(  \widetilde{P}\right)  =\left(  t_{1},t_{2}%
,\ldots,t_{n-k}\right)  $. Now, (\ref{sol.det.laplace-multi.0.short.a.kappa2})
shows that $\kappa\left(  k+i\right)  =t_{i}$ for every $i\in\left\{
1,2,\ldots,n-k\right\}  $. In other words, $\left(  \kappa\left(  k+1\right)
,\kappa\left(  k+2\right)  ,\ldots,\kappa\left(  n\right)  \right)  =\left(
t_{1},t_{2},\ldots,t_{n-k}\right)  $. Comparing this with $w\left(
\widetilde{P}\right)  =\left(  t_{1},t_{2},\ldots,t_{n-k}\right)  $, we obtain
$w\left(  \widetilde{P}\right)  =\left(  \kappa\left(  k+1\right)
,\kappa\left(  k+2\right)  ,\ldots,\kappa\left(  n\right)  \right)  $. This
proves (\ref{sol.det.laplace-multi.0.short.a.wPs=}).}.
\end{vershort}

\begin{verlong}
We have%
\begin{equation}
w\left(  R\right)  =\left(  \kappa\left(  1\right)  ,\kappa\left(  2\right)
,\ldots,\kappa\left(  k\right)  \right)  \label{sol.det.laplace-multi.0.a.wR=}%
\end{equation}
\footnote{\textit{Proof of (\ref{sol.det.laplace-multi.0.a.wR=}):} Recall that
$w\left(  R\right)  $ is the list of all elements of $R$ in increasing order
(with no repetitions) (by the definition of $w\left(  R\right)  $). Thus,%
\begin{align}
w\left(  R\right)   &  =\left(  \text{the list of all elements of }R\text{ in
increasing order (with no repetitions)}\right) \nonumber\\
&  =\left(  r_{1},r_{2},\ldots,r_{k}\right)
\label{sol.det.laplace-multi.0.a.wR=.pf.1}%
\end{align}
(since $\left(  r_{1},r_{2},\ldots,r_{k}\right)  $ is the list of all elements
of $R$ in increasing order (with no repetitions)). Now,
(\ref{sol.det.laplace-multi.0.a.kappa1}) shows that $\kappa\left(  i\right)
=r_{i}$ for every $i\in\left\{  1,2,\ldots,k\right\}  $. In other words,
$\left(  \kappa\left(  1\right)  ,\kappa\left(  2\right)  ,\ldots
,\kappa\left(  k\right)  \right)  =\left(  r_{1},r_{2},\ldots,r_{k}\right)  $.
Comparing this with (\ref{sol.det.laplace-multi.0.a.wR=.pf.1}), we obtain
$w\left(  R\right)  =\left(  \kappa\left(  1\right)  ,\kappa\left(  2\right)
,\ldots,\kappa\left(  k\right)  \right)  $. This proves
(\ref{sol.det.laplace-multi.0.a.wR=}).} and%
\begin{equation}
w\left(  \widetilde{P}\right)  =\left(  \kappa\left(  k+1\right)
,\kappa\left(  k+2\right)  ,\ldots,\kappa\left(  n\right)  \right)
\label{sol.det.laplace-multi.0.a.wPs=}%
\end{equation}
\footnote{\textit{Proof of (\ref{sol.det.laplace-multi.0.a.wPs=}):} Recall
that $w\left(  \widetilde{P}\right)  $ is the list of all elements of
$\widetilde{P}$ in increasing order (with no repetitions) (by the definition
of $w\left(  \widetilde{P}\right)  $). Thus,%
\begin{align}
w\left(  \widetilde{P}\right)   &  =\left(  \text{the list of all elements of
}\widetilde{P}\text{ in increasing order (with no repetitions)}\right)
\nonumber\\
&  =\left(  t_{1},t_{2},\ldots,t_{n-k}\right)
\label{sol.det.laplace-multi.0.a.wPs=.pf.1}%
\end{align}
(since $\left(  t_{1},t_{2},\ldots,t_{n-k}\right)  $ is the list of all
elements of $\widetilde{P}$ in increasing order (with no repetitions)). Now,
(\ref{sol.det.laplace-multi.0.a.kappa2}) shows that $\kappa\left(  k+i\right)
=t_{i}$ for every $i\in\left\{  1,2,\ldots,n-k\right\}  $. In other words,
$\left(  \kappa\left(  k+1\right)  ,\kappa\left(  k+2\right)  ,\ldots
,\kappa\left(  k+\left(  n-k\right)  \right)  \right)  =\left(  t_{1}%
,t_{2},\ldots,t_{n-k}\right)  $. Comparing this with
(\ref{sol.det.laplace-multi.0.a.wPs=.pf.1}), we obtain
\begin{align*}
w\left(  \widetilde{P}\right)   &  =\left(  \kappa\left(  k+1\right)
,\kappa\left(  k+2\right)  ,\ldots,\kappa\left(  k+\left(  n-k\right)
\right)  \right) \\
&  =\left(  \kappa\left(  k+1\right)  ,\kappa\left(  k+2\right)
,\ldots,\kappa\left(  n\right)  \right)  \ \ \ \ \ \ \ \ \ \ \left(
\text{since }k+\left(  n-k\right)  =n\right)  .
\end{align*}
This proves (\ref{sol.det.laplace-multi.0.a.wPs=}).}.
\end{verlong}

\begin{vershort}
Using $\widetilde{P}\cap R\neq\varnothing$, we can easily see that
$\kappa\notin S_{n}$\ \ \ \ \footnote{\textit{Proof.} Assume the contrary (for
the sake of contradiction). Thus, $\kappa\in S_{n}$. Hence, the inverse
permutation $\kappa^{-1}\in S_{n}$ is well-defined.
\par
Recall that $\widetilde{P}\cap R\neq\varnothing$. In other words, there exists
a $\rho\in\widetilde{P}\cap R$. Consider this $\rho$.
\par
Recall that $\left(  r_{1},r_{2},\ldots,r_{k}\right)  $ is a list of all
elements of $R$. Hence, $\left\{  r_{1},r_{2},\ldots,r_{k}\right\}  =R$. Now,
$\rho\in\widetilde{P}\cap R\subseteq R=\left\{  r_{1},r_{2},\ldots
,r_{k}\right\}  $. In other words, there exists some $i\in\left\{
1,2,\ldots,k\right\}  $ such that $\rho=r_{i}$. Consider this $i$.
\par
Recall that $\left(  t_{1},t_{2},\ldots,t_{n-k}\right)  $ is a list of all
elements of $\widetilde{P}$. Hence, $\left\{  t_{1},t_{2},\ldots
,t_{n-k}\right\}  =\widetilde{P}$. Now, $\rho\in\widetilde{P}\cap
R\subseteq\widetilde{P}=\left\{  t_{1},t_{2},\ldots,t_{n-k}\right\}  $. In
other words, there exists some $j\in\left\{  1,2,\ldots,n-k\right\}  $ such
that $\rho=t_{j}$. Consider this $j$.
\par
The equality (\ref{sol.det.laplace-multi.0.short.a.kappa2}) (applied to $j$
instead of $i$) yields $\kappa\left(  k+j\right)  =t_{j}=\rho$ (since
$\rho=t_{j}$). But from (\ref{sol.det.laplace-multi.0.short.a.kappa1}), we
obtain $\kappa\left(  i\right)  =r_{i}=\rho$ (since $\rho=r_{i}$). Comparing
this with $\kappa\left(  k+j\right)  =\rho$, we obtain $\kappa\left(
i\right)  =\kappa\left(  k+j\right)  $. Now, $i=\kappa^{-1}\left(
\underbrace{\kappa\left(  i\right)  }_{=\kappa\left(  k+j\right)  }\right)
=\kappa^{-1}\left(  \kappa\left(  k+j\right)  \right)  =k+\underbrace{j}%
_{\substack{>0\\\text{(since }j\in\left\{  1,2,\ldots,n-k\right\}  \text{)}%
}}>k$. But $i\in\left\{  1,2,\ldots,k\right\}  $ and thus $i\leq k$. This
contradicts $i>k$. This contradiction proves that our assumption was wrong.
Qed.}.
\end{vershort}

\begin{verlong}
Furthermore, $\kappa\notin S_{n}$\ \ \ \ \footnote{\textit{Proof.} Assume the
contrary (for the sake of contradiction). Thus, $\kappa\in S_{n}$.
\par
We have $\kappa\in S_{n}$. In other words, $\kappa$ is a permutation of
$\left\{  1,2,\ldots,n\right\}  $ (since $S_{n}$ is the set of all
permutations of $\left\{  1,2,\ldots,n\right\}  $). In other words, $\kappa$
is a bijective map $\left\{  1,2,\ldots,n\right\}  \rightarrow\left\{
1,2,\ldots,n\right\}  $. Thus, the map $\kappa$ is bijective, and therefore
injective.
\par
Recall that $\widetilde{P}\cap R\neq\varnothing$. Thus, the set $\widetilde{P}%
\cap R$ is nonempty. In other words, there exists a $\rho\in\widetilde{P}\cap
R$. Consider this $\rho$.
\par
Recall that $\left(  r_{1},r_{2},\ldots,r_{k}\right)  $ is the list of all
elements of $R$ in increasing order (with no repetitions). Thus, $\left(
r_{1},r_{2},\ldots,r_{k}\right)  $ is a list of all elements of $R$. Hence,
$\left\{  r_{1},r_{2},\ldots,r_{k}\right\}  =R$. Now, $\rho\in\widetilde{P}%
\cap R\subseteq R=\left\{  r_{1},r_{2},\ldots,r_{k}\right\}  $. In other
words, there exists some $i\in\left\{  1,2,\ldots,k\right\}  $ such that
$\rho=r_{i}$. Consider this $i$.
\par
Recall that $\left(  t_{1},t_{2},\ldots,t_{n-k}\right)  $ is the list of all
elements of $\widetilde{P}$ in increasing order (with no repetitions). Thus,
$\left(  t_{1},t_{2},\ldots,t_{n-k}\right)  $ is a list of all elements of
$\widetilde{P}$. Hence, $\left\{  t_{1},t_{2},\ldots,t_{n-k}\right\}
=\widetilde{P}$. Now, $\rho\in\widetilde{P}\cap R\subseteq\widetilde{P}%
=\left\{  t_{1},t_{2},\ldots,t_{n-k}\right\}  $. In other words, there exists
some $j\in\left\{  1,2,\ldots,n-k\right\}  $ such that $\rho=t_{j}$. Consider
this $j$.
\par
We have $j\in\left\{  1,2,\ldots,n-k\right\}  $, so that%
\begin{align*}
k+j  &  \in\left\{  k+1,k+2,\ldots,k+\left(  n-k\right)  \right\}  =\left\{
k+1,k+2,\ldots,n\right\}  \ \ \ \ \ \ \ \ \ \ \left(  \text{since }k+\left(
n-k\right)  =n\right) \\
&  \subseteq\left\{  1,2,\ldots,n\right\}  \ \ \ \ \ \ \ \ \ \ \left(
\text{since }k+1\geq1\right) \\
&  =\left[  n\right]  .
\end{align*}
Hence, $\kappa\left(  k+j\right)  $ is well-defined. The equality
(\ref{sol.det.laplace-multi.0.a.kappa2}) (applied to $j$ instead of $i$)
yields $\kappa\left(  k+j\right)  =t_{j}=\rho$ (since $\rho=t_{j}$).
\par
But%
\begin{align*}
i  &  \in\left\{  1,2,\ldots,k\right\}  \subseteq\left\{  1,2,\ldots
,n\right\}  \ \ \ \ \ \ \ \ \ \ \left(  \text{since }k\leq n\right) \\
&  =\left[  n\right]  .
\end{align*}
Hence, $\kappa\left(  i\right)  $ is well-defined. From
(\ref{sol.det.laplace-multi.0.a.kappa1}), we obtain $\kappa\left(  i\right)
=r_{i}=\rho$ (since $\rho=r_{i}$). Comparing this with $\kappa\left(
k+j\right)  =\rho$, we obtain $\kappa\left(  i\right)  =\kappa\left(
k+j\right)  $.
\par
But the map $\kappa$ is injective. In other words, if $u$ and $v$ are two
elements of $\left[  n\right]  $ satisfying $\kappa\left(  u\right)
=\kappa\left(  v\right)  $, then $u=v$. Applying this to $u=i$ and $v=k+j$, we
obtain $i=k+j$ (since $\kappa\left(  i\right)  =\kappa\left(  k+j\right)  $).
But $j\in\left\{  1,2,\ldots,n-k\right\}  $, so that $j\geq1>0$ and thus
$k+\underbrace{j}_{>0}>k$. Now, $i=k+j>k$. But $i\in\left\{  1,2,\ldots
,k\right\}  $ and thus $i\leq k$. This contradicts $i>k$. This contradiction
proves that our assumption was wrong. Qed.}.
\end{verlong}

Write the $n\times n$-matrix $A$ in the form $A=\left(  a_{i,j}\right)
_{1\leq i\leq n,\ 1\leq j\leq n}$. Define an $n\times n$-matrix $A_{\kappa}$
by $A_{\kappa}=\left(  a_{\kappa\left(  i\right)  ,j}\right)  _{1\leq i\leq
n,\ 1\leq j\leq n}$. Then, Lemma \ref{lem.det.sigma} \textbf{(b)} (applied to
$A$, $a_{i,j}$ and $A_{\kappa}$ instead of $B$, $b_{i,j}$ and $B_{\kappa}$)
yields
\begin{equation}
\det\left(  A_{\kappa}\right)  =0 \label{sol.det.laplace-multi.0.a.detAk}%
\end{equation}
(since $\kappa\notin S_{n}$).

Now, every subset $Q$ of $\left\{  1,2,\ldots,n\right\}  $ satisfies%
\begin{equation}
\operatorname*{sub}\nolimits_{w\left(  R\right)  }^{w\left(  Q\right)
}A=\operatorname*{sub}\nolimits_{\left(  1,2,\ldots,k\right)  }^{w\left(
Q\right)  }\left(  A_{\kappa}\right)  \label{sol.det.laplace-multi.0.a.sub1}%
\end{equation}
\footnote{\textit{Proof of (\ref{sol.det.laplace-multi.0.a.sub1}):} Let $Q$ be
a subset of $\left\{  1,2,\ldots,n\right\}  $. Write the list $w\left(
Q\right)  $ in the form $w\left(  Q\right)  =\left(  q_{1},q_{2}%
,\ldots,q_{\ell}\right)  $ for some $\ell\in\mathbb{N}$.
\par
From $w\left(  R\right)  =\left(  \kappa\left(  1\right)  ,\kappa\left(
2\right)  ,\ldots,\kappa\left(  k\right)  \right)  $ and $w\left(  Q\right)
=\left(  q_{1},q_{2},\ldots,q_{\ell}\right)  $, we obtain%
\begin{align}
\operatorname*{sub}\nolimits_{w\left(  R\right)  }^{w\left(  Q\right)  }A  &
=\operatorname*{sub}\nolimits_{\left(  \kappa\left(  1\right)  ,\kappa\left(
2\right)  ,\ldots,\kappa\left(  k\right)  \right)  }^{\left(  q_{1}%
,q_{2},\ldots,q_{\ell}\right)  }A=\operatorname*{sub}\nolimits_{\kappa\left(
1\right)  ,\kappa\left(  2\right)  ,\ldots,\kappa\left(  k\right)  }%
^{q_{1},q_{2},\ldots,q_{\ell}}A=\left(  a_{\kappa\left(  x\right)  ,q_{y}%
}\right)  _{1\leq x\leq k,\ 1\leq y\leq\ell}%
\label{sol.det.laplace-multi.0.a.sub1.pf.1}\\
&  \ \ \ \ \ \ \ \ \ \ \left(  \text{by the definition of }\operatorname*{sub}%
\nolimits_{\kappa\left(  1\right)  ,\kappa\left(  2\right)  ,\ldots
,\kappa\left(  k\right)  }^{q_{1},q_{2},\ldots,q_{\ell}}A\text{, since
}A=\left(  a_{i,j}\right)  _{1\leq i\leq n,\ 1\leq j\leq n}\right)  .\nonumber
\end{align}
\par
On the other hand, from $w\left(  Q\right)  =\left(  q_{1},q_{2}%
,\ldots,q_{\ell}\right)  $, we obtain%
\[
\operatorname*{sub}\nolimits_{\left(  1,2,\ldots,k\right)  }^{w\left(
Q\right)  }\left(  A_{\kappa}\right)  =\operatorname*{sub}\nolimits_{\left(
1,2,\ldots,k\right)  }^{\left(  q_{1},q_{2},\ldots,q_{\ell}\right)  }\left(
A_{\kappa}\right)  =\operatorname*{sub}\nolimits_{1,2,\ldots,k}^{q_{1}%
,q_{2},\ldots,q_{\ell}}\left(  A_{\kappa}\right)  =\left(  a_{\kappa\left(
x\right)  ,q_{y}}\right)  _{1\leq x\leq k,\ 1\leq y\leq\ell}%
\]
(by the definition of $\operatorname*{sub}\nolimits_{1,2,\ldots,k}%
^{q_{1},q_{2},\ldots,q_{\ell}}\left(  A_{\kappa}\right)  $, since $A_{\kappa
}=\left(  a_{\kappa\left(  i\right)  ,j}\right)  _{1\leq i\leq n,\ 1\leq j\leq
n}$). Comparing this with (\ref{sol.det.laplace-multi.0.a.sub1.pf.1}), we
obtain $\operatorname*{sub}\nolimits_{w\left(  R\right)  }^{w\left(  Q\right)
}A=\operatorname*{sub}\nolimits_{\left(  1,2,\ldots,k\right)  }^{w\left(
Q\right)  }\left(  A_{\kappa}\right)  $. This proves
(\ref{sol.det.laplace-multi.0.a.sub1}).} and%
\begin{equation}
\operatorname*{sub}\nolimits_{w\left(  \widetilde{P}\right)  }^{w\left(
\widetilde{Q}\right)  }A=\operatorname*{sub}\nolimits_{\left(  k+1,k+2,\ldots
,n\right)  }^{w\left(  \widetilde{Q}\right)  }\left(  A_{\kappa}\right)
\label{sol.det.laplace-multi.0.a.sub2}%
\end{equation}
\footnote{\textit{Proof of (\ref{sol.det.laplace-multi.0.a.sub2}):} Let $Q$ be
a subset of $\left\{  1,2,\ldots,n\right\}  $. Write the list $w\left(
\widetilde{Q}\right)  $ in the form $w\left(  \widetilde{Q}\right)  =\left(
q_{1},q_{2},\ldots,q_{\ell}\right)  $ for some $\ell\in\mathbb{N}$.
\par
From $w\left(  \widetilde{P}\right)  =\left(  \kappa\left(  k+1\right)
,\kappa\left(  k+2\right)  ,\ldots,\kappa\left(  n\right)  \right)  $ and
$w\left(  \widetilde{Q}\right)  =\left(  q_{1},q_{2},\ldots,q_{\ell}\right)
$, we obtain%
\begin{align}
\operatorname*{sub}\nolimits_{w\left(  \widetilde{P}\right)  }^{w\left(
\widetilde{Q}\right)  }A  &  =\operatorname*{sub}\nolimits_{\left(
\kappa\left(  k+1\right)  ,\kappa\left(  k+2\right)  ,\ldots,\kappa\left(
n\right)  \right)  }^{\left(  q_{1},q_{2},\ldots,q_{\ell}\right)
}A=\operatorname*{sub}\nolimits_{\kappa\left(  k+1\right)  ,\kappa\left(
k+2\right)  ,\ldots,\kappa\left(  n\right)  }^{q_{1},q_{2},\ldots,q_{\ell}%
}A=\left(  a_{\kappa\left(  k+x\right)  ,q_{y}}\right)  _{1\leq x\leq
n-k,\ 1\leq y\leq\ell}\label{sol.det.laplace-multi.0.a.sub2.pf.1}\\
&  \ \ \ \ \ \ \ \ \ \ \left(  \text{by the definition of }\operatorname*{sub}%
\nolimits_{\kappa\left(  k+1\right)  ,\kappa\left(  k+2\right)  ,\ldots
,\kappa\left(  n\right)  }^{q_{1},q_{2},\ldots,q_{\ell}}A\text{, since
}A=\left(  a_{i,j}\right)  _{1\leq i\leq n,\ 1\leq j\leq n}\right)  .\nonumber
\end{align}
\par
On the other hand, from $w\left(  \widetilde{Q}\right)  =\left(  q_{1}%
,q_{2},\ldots,q_{\ell}\right)  $, we obtain%
\[
\operatorname*{sub}\nolimits_{\left(  k+1,k+2,\ldots,n\right)  }^{w\left(
\widetilde{Q}\right)  }\left(  A_{\kappa}\right)  =\operatorname*{sub}%
\nolimits_{\left(  k+1,k+2,\ldots,n\right)  }^{\left(  q_{1},q_{2}%
,\ldots,q_{\ell}\right)  }\left(  A_{\kappa}\right)  =\operatorname*{sub}%
\nolimits_{k+1,k+2,\ldots,n}^{q_{1},q_{2},\ldots,q_{\ell}}\left(  A_{\kappa
}\right)  =\left(  a_{\kappa\left(  k+x\right)  ,q_{y}}\right)  _{1\leq x\leq
n-k,\ 1\leq y\leq\ell}%
\]
(by the definition of $\operatorname*{sub}\nolimits_{k+1,k+2,\ldots,n}%
^{q_{1},q_{2},\ldots,q_{\ell}}\left(  A_{\kappa}\right)  $, since $A_{\kappa
}=\left(  a_{\kappa\left(  i\right)  ,j}\right)  _{1\leq i\leq n,\ 1\leq j\leq
n}$). Comparing this with (\ref{sol.det.laplace-multi.0.a.sub2.pf.1}), we
obtain $\operatorname*{sub}\nolimits_{w\left(  \widetilde{P}\right)
}^{w\left(  \widetilde{Q}\right)  }A=\operatorname*{sub}\nolimits_{\left(
k+1,k+2,\ldots,n\right)  }^{w\left(  \widetilde{Q}\right)  }\left(  A_{\kappa
}\right)  $. This proves (\ref{sol.det.laplace-multi.0.a.sub2}).}.

\begin{vershort}
Define a subset $P^{\prime}$ of $\left\{  1,2,\ldots,n\right\}  $ by
$P^{\prime}=\left\{  1,2,\ldots,k\right\}  $. Then, $\left\vert P^{\prime
}\right\vert =k=\left\vert P\right\vert $.
\end{vershort}

\begin{verlong}
Now, we have $k\leq n$ and thus $\left\{  1,2,\ldots,k\right\}  \subseteq
\left\{  1,2,\ldots,n\right\}  $. In other words, $\left\{  1,2,\ldots
,k\right\}  $ is a subset of $\left\{  1,2,\ldots,n\right\}  $. Thus, we can
define a subset $P^{\prime}$ of $\left\{  1,2,\ldots,n\right\}  $ by
$P^{\prime}=\left\{  1,2,\ldots,k\right\}  $. Consider this $P^{\prime}$. From
$P^{\prime}=\left\{  1,2,\ldots,k\right\}  $, we obtain $\left\vert P^{\prime
}\right\vert =\left\vert \left\{  1,2,\ldots,k\right\}  \right\vert
=k=\left\vert P\right\vert $.
\end{verlong}

\begin{vershort}
We have $w\left(  P^{\prime}\right)  =\left(  1,2,\ldots,k\right)
$\ \ \ \ \footnote{\textit{Proof.} The definition of $w\left(  P^{\prime
}\right)  $ yields%
\begin{align*}
w\left(  P^{\prime}\right)   &  =\left(  \text{the list of all elements of
}\underbrace{P^{\prime}}_{=\left\{  1,2,\ldots,k\right\}  }\text{ in
increasing order (with no repetitions)}\right) \\
&  =\left(  \text{the list of all elements of }\left\{  1,2,\ldots,k\right\}
\text{ in increasing order (with no repetitions)}\right) \\
&  =\left(  1,2,\ldots,k\right)  .
\end{align*}
Qed.} and $w\left(  \widetilde{P^{\prime}}\right)  =\left(  k+1,k+2,\ldots
,n\right)  $\ \ \ \ \footnote{\textit{Proof.} We have $P^{\prime}=\left\{
1,2,\ldots,k\right\}  $. Now,
\begin{align*}
\widetilde{P^{\prime}}  &  =\left\{  1,2,\ldots,n\right\}  \setminus
\underbrace{P^{\prime}}_{=\left\{  1,2,\ldots,k\right\}  }%
\ \ \ \ \ \ \ \ \ \ \left(  \text{by the definition of }\widetilde{P^{\prime}%
}\right) \\
&  =\left\{  1,2,\ldots,n\right\}  \setminus\left\{  1,2,\ldots,k\right\}
=\left\{  k+1,k+2,\ldots,n\right\}  \ \ \ \ \ \ \ \ \ \ \left(  \text{since
}k\in\left\{  0,1,\ldots,n\right\}  \right)  .
\end{align*}
\par
Now, the definition of $w\left(  \widetilde{P^{\prime}}\right)  $ yields%
\begin{align*}
&  w\left(  \widetilde{P^{\prime}}\right) \\
&  =\left(  \text{the list of all elements of }%
\underbrace{\widetilde{P^{\prime}}}_{=\left\{  k+1,k+2,\ldots,n\right\}
}\text{ in increasing order (with no repetitions)}\right) \\
&  =\left(  \text{the list of all elements of }\left\{  k+1,k+2,\ldots
,n\right\}  \text{ in increasing order (with no repetitions)}\right) \\
&  =\left(  k+1,k+2,\ldots,n\right)  .
\end{align*}
Qed.}.
\end{vershort}

\begin{verlong}
We have $w\left(  P^{\prime}\right)  =\left(  1,2,\ldots,k\right)
$\ \ \ \ \footnote{\textit{Proof.} Recall that $w\left(  P^{\prime}\right)  $
is the list of all elements of $P^{\prime}$ in increasing order (with no
repetitions) (by the definition of $w\left(  P^{\prime}\right)  $). Thus,%
\begin{align*}
w\left(  P^{\prime}\right)   &  =\left(  \text{the list of all elements of
}\underbrace{P^{\prime}}_{=\left\{  1,2,\ldots,k\right\}  }\text{ in
increasing order (with no repetitions)}\right) \\
&  =\left(  \text{the list of all elements of }\left\{  1,2,\ldots,k\right\}
\text{ in increasing order (with no repetitions)}\right) \\
&  =\left(  1,2,\ldots,k\right)  .
\end{align*}
Qed.} and $w\left(  \widetilde{P^{\prime}}\right)  =\left(  k+1,k+2,\ldots
,n\right)  $\ \ \ \ \footnote{\textit{Proof.} We have $P^{\prime}=\left\{
1,2,\ldots,k\right\}  $. Now,
\begin{align*}
\widetilde{P^{\prime}}  &  =\left\{  1,2,\ldots,n\right\}  \setminus
\underbrace{P^{\prime}}_{=\left\{  1,2,\ldots,k\right\}  }%
\ \ \ \ \ \ \ \ \ \ \left(  \text{by the definition of }\widetilde{P^{\prime}%
}\right) \\
&  =\left\{  1,2,\ldots,n\right\}  \setminus\left\{  1,2,\ldots,k\right\}
=\left\{  k+1,k+2,\ldots,n\right\}  \ \ \ \ \ \ \ \ \ \ \left(  \text{since
}k\in\left\{  0,1,\ldots,n\right\}  \right)  .
\end{align*}
\par
Recall that $w\left(  \widetilde{P^{\prime}}\right)  $ is the list of all
elements of $\widetilde{P^{\prime}}$ in increasing order (with no repetitions)
(by the definition of $w\left(  \widetilde{P^{\prime}}\right)  $). Thus,%
\begin{align*}
&  w\left(  \widetilde{P^{\prime}}\right) \\
&  =\left(  \text{the list of all elements of }%
\underbrace{\widetilde{P^{\prime}}}_{=\left\{  k+1,k+2,\ldots,n\right\}
}\text{ in increasing order (with no repetitions)}\right) \\
&  =\left(  \text{the list of all elements of }\left\{  k+1,k+2,\ldots
,n\right\}  \text{ in increasing order (with no repetitions)}\right) \\
&  =\left(  k+1,k+2,\ldots,n\right)  .
\end{align*}
Qed.}.
\end{verlong}

Now, Theorem \ref{thm.det.laplace-multi} \textbf{(a)} (applied to $A_{\kappa}$
and $P^{\prime}$ instead of $A$ and $P$) yields%
\begin{align*}
&  \det\left(  A_{\kappa}\right) \\
&  =\underbrace{\sum_{\substack{Q\subseteq\left\{  1,2,\ldots,n\right\}
;\\\left\vert Q\right\vert =\left\vert P^{\prime}\right\vert }}}%
_{\substack{=\sum_{\substack{Q\subseteq\left\{  1,2,\ldots,n\right\}
;\\\left\vert Q\right\vert =\left\vert P\right\vert }}\\\text{(since
}\left\vert P^{\prime}\right\vert =\left\vert P\right\vert \text{)}}}\left(
-1\right)  ^{\sum P^{\prime}+\sum Q}\det\left(
\underbrace{\operatorname*{sub}\nolimits_{w\left(  P^{\prime}\right)
}^{w\left(  Q\right)  }\left(  A_{\kappa}\right)  }%
_{\substack{=\operatorname*{sub}\nolimits_{\left(  1,2,\ldots,k\right)
}^{w\left(  Q\right)  }\left(  A_{\kappa}\right)  \\\text{(since }w\left(
P^{\prime}\right)  =\left(  1,2,\ldots,k\right)  \text{)}}}\right)
\det\left(  \underbrace{\operatorname*{sub}\nolimits_{w\left(
\widetilde{P^{\prime}}\right)  }^{w\left(  \widetilde{Q}\right)  }\left(
A_{\kappa}\right)  }_{\substack{=\operatorname*{sub}\nolimits_{\left(
k+1,k+2,\ldots,n\right)  }^{w\left(  \widetilde{Q}\right)  }\left(  A_{\kappa
}\right)  \\\text{(since }w\left(  \widetilde{P^{\prime}}\right)  =\left(
k+1,k+2,\ldots,n\right)  \text{)}}}\right) \\
&  =\sum_{\substack{Q\subseteq\left\{  1,2,\ldots,n\right\}  ;\\\left\vert
Q\right\vert =\left\vert P\right\vert }}\left(  -1\right)  ^{\sum P^{\prime
}+\sum Q}\det\left(  \underbrace{\operatorname*{sub}\nolimits_{\left(
1,2,\ldots,k\right)  }^{w\left(  Q\right)  }\left(  A_{\kappa}\right)
}_{\substack{=\operatorname*{sub}\nolimits_{w\left(  R\right)  }^{w\left(
Q\right)  }A\\\text{(by (\ref{sol.det.laplace-multi.0.a.sub1}))}}}\right)
\det\left(  \underbrace{\operatorname*{sub}\nolimits_{\left(  k+1,k+2,\ldots
,n\right)  }^{w\left(  \widetilde{Q}\right)  }\left(  A_{\kappa}\right)
}_{\substack{=\operatorname*{sub}\nolimits_{w\left(  \widetilde{P}\right)
}^{w\left(  \widetilde{Q}\right)  }A\\\text{(by
(\ref{sol.det.laplace-multi.0.a.sub2}))}}}\right) \\
&  =\sum_{\substack{Q\subseteq\left\{  1,2,\ldots,n\right\}  ;\\\left\vert
Q\right\vert =\left\vert P\right\vert }}\left(  -1\right)  ^{\sum P^{\prime
}+\sum Q}\det\left(  \operatorname*{sub}\nolimits_{w\left(  R\right)
}^{w\left(  Q\right)  }A\right)  \det\left(  \operatorname*{sub}%
\nolimits_{w\left(  \widetilde{P}\right)  }^{w\left(  \widetilde{Q}\right)
}A\right)  .
\end{align*}
Comparing this with (\ref{sol.det.laplace-multi.0.a.detAk}), we obtain%
\[
0=\sum_{\substack{Q\subseteq\left\{  1,2,\ldots,n\right\}  ;\\\left\vert
Q\right\vert =\left\vert P\right\vert }}\left(  -1\right)  ^{\sum P^{\prime
}+\sum Q}\det\left(  \operatorname*{sub}\nolimits_{w\left(  R\right)
}^{w\left(  Q\right)  }A\right)  \det\left(  \operatorname*{sub}%
\nolimits_{w\left(  \widetilde{P}\right)  }^{w\left(  \widetilde{Q}\right)
}A\right)  .
\]
Multiplying both sides of this equality by $\left(  -1\right)  ^{\sum P-\sum
P^{\prime}}$, we obtain%
\begin{align*}
0  &  =\left(  -1\right)  ^{\sum P-\sum P^{\prime}}\cdot\sum
_{\substack{Q\subseteq\left\{  1,2,\ldots,n\right\}  ;\\\left\vert
Q\right\vert =\left\vert P\right\vert }}\left(  -1\right)  ^{\sum P^{\prime
}+\sum Q}\det\left(  \operatorname*{sub}\nolimits_{w\left(  R\right)
}^{w\left(  Q\right)  }A\right)  \det\left(  \operatorname*{sub}%
\nolimits_{w\left(  \widetilde{P}\right)  }^{w\left(  \widetilde{Q}\right)
}A\right) \\
&  =\sum_{\substack{Q\subseteq\left\{  1,2,\ldots,n\right\}  ;\\\left\vert
Q\right\vert =\left\vert P\right\vert }}\underbrace{\left(  -1\right)  ^{\sum
P-\sum P^{\prime}}\left(  -1\right)  ^{\sum P^{\prime}+\sum Q}}%
_{\substack{=\left(  -1\right)  ^{\left(  \sum P-\sum P^{\prime}\right)
+\left(  \sum P^{\prime}+\sum Q\right)  }\\=\left(  -1\right)  ^{\sum P+\sum
Q}\\\text{(since }\left(  \sum P-\sum P^{\prime}\right)  +\left(  \sum
P^{\prime}+\sum Q\right)  =\sum P+\sum Q\text{)}}}\det\left(
\operatorname*{sub}\nolimits_{w\left(  R\right)  }^{w\left(  Q\right)
}A\right)  \det\left(  \operatorname*{sub}\nolimits_{w\left(  \widetilde{P}%
\right)  }^{w\left(  \widetilde{Q}\right)  }A\right) \\
&  =\sum_{\substack{Q\subseteq\left\{  1,2,\ldots,n\right\}  ;\\\left\vert
Q\right\vert =\left\vert P\right\vert }}\left(  -1\right)  ^{\sum P+\sum
Q}\det\left(  \operatorname*{sub}\nolimits_{w\left(  R\right)  }^{w\left(
Q\right)  }A\right)  \det\left(  \operatorname*{sub}\nolimits_{w\left(
\widetilde{P}\right)  }^{w\left(  \widetilde{Q}\right)  }A\right)  .
\end{align*}
This solves Exercise \ref{exe.det.laplace-multi.0} \textbf{(a)}.

\textbf{(b)} Let $P$ be a subset of $\left\{  1,2,\ldots,n\right\}  $
satisfying $\left\vert P\right\vert =\left\vert R\right\vert $ and $P\neq R$.
Every subset $Q$ of $\left\{  1,2,\ldots,n\right\}  $ satisfying $\left\vert
Q\right\vert =\left\vert P\right\vert $ satisfies%
\begin{align}
&  \det\left(  \operatorname*{sub}\nolimits_{w\left(  R\right)  }^{w\left(
Q\right)  }\left(  A^{T}\right)  \right)  \det\left(  \operatorname*{sub}%
\nolimits_{w\left(  \widetilde{P}\right)  }^{w\left(  \widetilde{Q}\right)
}\left(  A^{T}\right)  \right) \nonumber\\
&  =\det\left(  \operatorname*{sub}\nolimits_{w\left(  Q\right)  }^{w\left(
R\right)  }A\right)  \det\left(  \operatorname*{sub}\nolimits_{w\left(
\widetilde{Q}\right)  }^{w\left(  \widetilde{P}\right)  }A\right)
\label{sol.det.laplace-multi.0.b.1}%
\end{align}
\footnote{\textit{Proof of (\ref{sol.det.laplace-multi.0.b.1}):} Let $Q$ be a
subset of $\left\{  1,2,\ldots,n\right\}  $ satisfying $\left\vert
Q\right\vert =\left\vert P\right\vert $. Then, $\left\vert Q\right\vert
=\left\vert P\right\vert =\left\vert R\right\vert $. Thus, $\left\vert
R\right\vert =\left\vert Q\right\vert $. Hence, Corollary
\ref{cor.sol.det.laplace-multi.0.detsubAT} (applied to $m=n$, $U=R$ and $V=Q$)
yields
\begin{equation}
\det\left(  \operatorname*{sub}\nolimits_{w\left(  R\right)  }^{w\left(
Q\right)  }\left(  A^{T}\right)  \right)  =\det\left(  \operatorname*{sub}%
\nolimits_{w\left(  Q\right)  }^{w\left(  R\right)  }A\right)  .
\label{sol.det.laplace-multi.0.b.1.pf.1a}%
\end{equation}
\par
On the other hand, the definition of $\widetilde{P}$ yields $\widetilde{P}%
=\left\{  1,2,\ldots,n\right\}  \setminus P\subseteq\left\{  1,2,\ldots
,n\right\}  $. Hence, $\widetilde{P}$ is a subset of $\left\{  1,2,\ldots
,n\right\}  $. The same argument (applied to $Q$ instead of $P$) shows that
$\widetilde{Q}$ is a subset of $\left\{  1,2,\ldots,n\right\}  $. Furthermore,
$\widetilde{P}=\left\{  1,2,\ldots,n\right\}  \setminus P$ and thus
\begin{align*}
\left\vert \widetilde{P}\right\vert  &  =\left\vert \left\{  1,2,\ldots
,n\right\}  \setminus P\right\vert =\underbrace{\left\vert \left\{
1,2,\ldots,n\right\}  \right\vert }_{=n}-\left\vert P\right\vert
\ \ \ \ \ \ \ \ \ \ \left(  \text{since }P\text{ is a subset of }\left\{
1,2,\ldots,n\right\}  \right) \\
&  =n-\left\vert P\right\vert .
\end{align*}
The same argument (applied to $Q$ instead of $P$) shows that $\left\vert
\widetilde{Q}\right\vert =n-\left\vert Q\right\vert $. Hence, $\left\vert
\widetilde{Q}\right\vert =n-\underbrace{\left\vert Q\right\vert }_{=\left\vert
P\right\vert }=n-\left\vert P\right\vert $. Comparing this with $\left\vert
\widetilde{P}\right\vert =n-\left\vert P\right\vert $, we obtain $\left\vert
\widetilde{P}\right\vert =\left\vert \widetilde{Q}\right\vert $. Thus,
Corollary \ref{cor.sol.det.laplace-multi.0.detsubAT} (applied to $m=n$,
$U=\widetilde{P}$ and $V=\widetilde{Q}$) yields
\begin{equation}
\det\left(  \operatorname*{sub}\nolimits_{w\left(  \widetilde{P}\right)
}^{w\left(  \widetilde{Q}\right)  }\left(  A^{T}\right)  \right)  =\det\left(
\operatorname*{sub}\nolimits_{w\left(  \widetilde{Q}\right)  }^{w\left(
\widetilde{P}\right)  }A\right)  . \label{sol.det.laplace-multi.0.b.1.pf.2a}%
\end{equation}
Multiplying the equality (\ref{sol.det.laplace-multi.0.b.1.pf.1a}) with the
equality (\ref{sol.det.laplace-multi.0.b.1.pf.2a}), we obtain%
\[
\det\left(  \operatorname*{sub}\nolimits_{w\left(  R\right)  }^{w\left(
Q\right)  }\left(  A^{T}\right)  \right)  \det\left(  \operatorname*{sub}%
\nolimits_{w\left(  \widetilde{P}\right)  }^{w\left(  \widetilde{Q}\right)
}\left(  A^{T}\right)  \right)  =\det\left(  \operatorname*{sub}%
\nolimits_{w\left(  Q\right)  }^{w\left(  R\right)  }A\right)  \det\left(
\operatorname*{sub}\nolimits_{w\left(  \widetilde{Q}\right)  }^{w\left(
\widetilde{P}\right)  }A\right)  .
\]
Thus, (\ref{sol.det.laplace-multi.0.b.1}) is proven.}.

Exercise \ref{exe.det.laplace-multi.0} \textbf{(a)} (applied to $A^{T}$
instead of $A$) yields%
\begin{align}
0  &  =\sum_{\substack{Q\subseteq\left\{  1,2,\ldots,n\right\}  ;\\\left\vert
Q\right\vert =\left\vert P\right\vert }}\underbrace{\left(  -1\right)  ^{\sum
P+\sum Q}}_{\substack{=\left(  -1\right)  ^{\sum Q+\sum P}\\\text{(since }\sum
P+\sum Q=\sum Q+\sum P\text{)}}}\underbrace{\det\left(  \operatorname*{sub}%
\nolimits_{w\left(  R\right)  }^{w\left(  Q\right)  }\left(  A^{T}\right)
\right)  \det\left(  \operatorname*{sub}\nolimits_{w\left(  \widetilde{P}%
\right)  }^{w\left(  \widetilde{Q}\right)  }\left(  A^{T}\right)  \right)
}_{\substack{=\det\left(  \operatorname*{sub}\nolimits_{w\left(  Q\right)
}^{w\left(  R\right)  }A\right)  \det\left(  \operatorname*{sub}%
\nolimits_{w\left(  \widetilde{Q}\right)  }^{w\left(  \widetilde{P}\right)
}A\right)  \\\text{(by (\ref{sol.det.laplace-multi.0.b.1}))}}}\nonumber\\
&  =\sum_{\substack{Q\subseteq\left\{  1,2,\ldots,n\right\}  ;\\\left\vert
Q\right\vert =\left\vert P\right\vert }}\left(  -1\right)  ^{\sum Q+\sum
P}\det\left(  \operatorname*{sub}\nolimits_{w\left(  Q\right)  }^{w\left(
R\right)  }A\right)  \det\left(  \operatorname*{sub}\nolimits_{w\left(
\widetilde{Q}\right)  }^{w\left(  \widetilde{P}\right)  }A\right) \nonumber\\
&  =\sum_{\substack{G\subseteq\left\{  1,2,\ldots,n\right\}  ;\\\left\vert
G\right\vert =\left\vert P\right\vert }}\left(  -1\right)  ^{\sum G+\sum
P}\det\left(  \operatorname*{sub}\nolimits_{w\left(  G\right)  }^{w\left(
R\right)  }A\right)  \det\left(  \operatorname*{sub}\nolimits_{w\left(
\widetilde{G}\right)  }^{w\left(  \widetilde{P}\right)  }A\right)
\label{sol.det.laplace-multi.0.b.5}%
\end{align}
(here, we have renamed the summation index $Q$ as $G$).

Now, forget that we fixed $P$. We thus have proven that
(\ref{sol.det.laplace-multi.0.b.5}) holds for every subset $P$ of $\left\{
1,2,\ldots,n\right\}  $ satisfying $\left\vert P\right\vert =\left\vert
R\right\vert $ and $P\neq R$.

Now, let $Q$ be a subset of $\left\{  1,2,\ldots,n\right\}  $ satisfying
$\left\vert Q\right\vert =\left\vert R\right\vert $ and $Q\neq R$. Then, we
can apply (\ref{sol.det.laplace-multi.0.b.5}) to $P=Q$. We thus obtain%
\begin{align*}
0  &  =\sum_{\substack{G\subseteq\left\{  1,2,\ldots,n\right\}  ;\\\left\vert
G\right\vert =\left\vert Q\right\vert }}\left(  -1\right)  ^{\sum G+\sum
Q}\det\left(  \operatorname*{sub}\nolimits_{w\left(  G\right)  }^{w\left(
R\right)  }A\right)  \det\left(  \operatorname*{sub}\nolimits_{w\left(
\widetilde{G}\right)  }^{w\left(  \widetilde{Q}\right)  }A\right) \\
&  =\sum_{\substack{P\subseteq\left\{  1,2,\ldots,n\right\}  ;\\\left\vert
P\right\vert =\left\vert Q\right\vert }}\left(  -1\right)  ^{\sum P+\sum
Q}\det\left(  \operatorname*{sub}\nolimits_{w\left(  P\right)  }^{w\left(
Q\right)  }A\right)  \det\left(  \operatorname*{sub}\nolimits_{w\left(
\widetilde{P}\right)  }^{w\left(  \widetilde{Q}\right)  }A\right)
\end{align*}
(here, we have renamed the summation index $G$ as $P$). This solves Exercise
\ref{exe.det.laplace-multi.0} \textbf{(b)}.
\end{proof}

\subsection{Solution to Exercise \ref{exe.det.laplace-multi.0r}}

Throughout this section, we shall use the notations introduced in Definition
\ref{def.submatrix} and in Definition \ref{def.sect.laplace.notations}. Also,
whenever $m$ is an integer, we shall use the notation $\left[  m\right]  $ for
the set $\left\{  1,2,\ldots,m\right\}  $.

\begin{proof}
[Solution to Exercise \ref{exe.det.laplace-multi.0r}.]\textbf{(a)} Let
$A\in\mathbb{K}^{m\times n}$.

\begin{vershort}
Define a $p\in\mathbb{N}$ by $p=\left\vert J\right\vert $. From $\left\vert
J\right\vert +\left\vert K\right\vert =n$, we obtain $\left\vert K\right\vert
=n-\underbrace{\left\vert J\right\vert }_{=p}=n-p$. Hence, $n-p=\left\vert
K\right\vert \geq0$, so that $n\geq p$ and therefore $p\leq n$. Thus,
$p\in\left\{  0,1,\ldots,n\right\}  $ and $\left\{  1,2,\ldots,p\right\}
\subseteq\left\{  1,2,\ldots,n\right\}  $ and $\left\{  p+1,p+2,\ldots
,n\right\}  \subseteq\left\{  1,2,\ldots,n\right\}  $.
\end{vershort}

\begin{verlong}
We have $\left[  m\right]  =\left\{  1,2,\ldots,m\right\}  $ (by the
definition of $\left[  m\right]  $) and $\left[  n\right]  =\left\{
1,2,\ldots,n\right\}  $ (by the definition of $\left[  n\right]  $).

The set $J$ is finite (since $J$ is a subset of the finite set $\left\{
1,2,\ldots,m\right\}  $). Define a $p\in\mathbb{N}$ by $p=\left\vert
J\right\vert $. (This is well-defined, since the set $J$ is finite.) From
$p\in\mathbb{N}$, we obtain $p\geq0$.

We have $\left\vert J\right\vert +\left\vert K\right\vert =n$, so that
$\left\vert K\right\vert =n-\underbrace{\left\vert J\right\vert }_{=p}=n-p$.
Hence, $n-p=\left\vert K\right\vert \geq0$, so that $n\geq p$ and therefore
$p\leq n$. Thus, $\left\{  1,2,\ldots,p\right\}  \subseteq\left\{
1,2,\ldots,n\right\}  $ (since $p\geq0$). Also, from $\underbrace{p}_{\geq
0}+1\geq1$, we obtain $\left\{  p+1,p+2,\ldots,n\right\}  \subseteq\left\{
1,2,\ldots,n\right\}  $. Also, from $0\leq p\leq n$, we obtain $p\in\left\{
0,1,\ldots,n\right\}  $.
\end{verlong}

Let $\left(  j_{1},j_{2},\ldots,j_{p}\right)  $ be the list of all elements of
$J$ in increasing order (with no repetitions). (This is well-defined (by
Definition \ref{def.ind.inclist}), because $\left\vert J\right\vert =p$.)

Let $\left(  k_{1},k_{2},\ldots,k_{n-p}\right)  $ be the list of all elements
of $K$ in increasing order (with no repetitions). (This is well-defined,
because $\left\vert K\right\vert =n-p$.)

\begin{vershort}
The $p$ elements $j_{1},j_{2},\ldots,j_{p}$ all belong to $J$ (since $\left(
j_{1},j_{2},\ldots,j_{p}\right)  $ is a list of elements of $J$), and thus to
$\left[  m\right]  $ (since $J$ is a subset of $\left\{  1,2,\ldots,m\right\}
=\left[  m\right]  $). Similarly, the $n-p$ elements $k_{1},k_{2}%
,\ldots,k_{n-p}$ belong to $\left[  m\right]  $. Combining the previous two
sentences, we conclude that the $p+\left(  n-p\right)  $ elements $j_{1}%
,j_{2},\ldots,j_{p},k_{1},k_{2},\ldots,k_{n-p}$ belong to $\left[  m\right]
$. Hence, $\left(  j_{1},j_{2},\ldots,j_{p},k_{1},k_{2},\ldots,k_{n-p}\right)
\in\left[  m\right]  ^{p+\left(  n-p\right)  }=\left[  m\right]  ^{n}$ (since
$p+\left(  n-p\right)  =n$).

Thus, we can define an $n$-tuple $\left(  \gamma_{1},\gamma_{2},\ldots
,\gamma_{n}\right)  \in\left[  m\right]  ^{n}$ by
\begin{equation}
\left(  \gamma_{1},\gamma_{2},\ldots,\gamma_{n}\right)  =\left(  j_{1}%
,j_{2},\ldots,j_{p},k_{1},k_{2},\ldots,k_{n-p}\right)  .
\label{sol.det.laplace-multi.0r.short.1}%
\end{equation}
Consider this $\left(  \gamma_{1},\gamma_{2},\ldots,\gamma_{n}\right)  $. From
(\ref{sol.det.laplace-multi.0r.short.1}), we obtain%
\begin{equation}
\gamma_{i}=j_{i}\ \ \ \ \ \ \ \ \ \ \text{for each }i\in\left\{
1,2,\ldots,p\right\}  \label{sol.det.laplace-multi.0r.short.2a}%
\end{equation}
and%
\begin{equation}
\gamma_{p+i}=k_{i}\ \ \ \ \ \ \ \ \ \ \text{for each }i\in\left\{
1,2,\ldots,n-p\right\}  . \label{sol.det.laplace-multi.0r.short.2b}%
\end{equation}

We have $\gamma_{i}\in\left[  m\right]  $ for each $i\in\left[  n\right]  $
(since $\left(  \gamma_{1},\gamma_{2},\ldots,\gamma_{n}\right)  \in\left[
m\right]  ^{n}$). Hence, we can define a map $\gamma:\left[  n\right]
\rightarrow\left[  m\right]  $ by%
\[
\left(  \gamma\left(  i\right)  =\gamma_{i}\text{ for every }i\in\left[
n\right]  \right)  .
\]
Consider this $\gamma$.
\end{vershort}

\begin{verlong}
We know that $\left(  j_{1},j_{2},\ldots,j_{p}\right)  $ is a list of elements
of $J$ (since $\left(  j_{1},j_{2},\ldots,j_{p}\right)  $ is the list of all
elements of $J$ in increasing order (with no repetitions)). Hence, the $p$
elements $j_{1},j_{2},\ldots,j_{p}$ belong to $J$, and thus belong to $\left[
m\right]  $ as well (since $J\subseteq\left\{  1,2,\ldots,m\right\}  =\left[
m\right]  $).

We know that $\left(  k_{1},k_{2},\ldots,k_{n-p}\right)  $ is a list of
elements of $K$ (since $\left(  k_{1},k_{2},\ldots,k_{n-p}\right)  $ is the
list of all elements of $K$ in increasing order (with no repetitions)). Hence,
the $n-p$ elements $k_{1},k_{2},\ldots,k_{n-p}$ belong to $K$, and thus belong
to $\left[  m\right]  $ as well (since $K\subseteq\left\{  1,2,\ldots
,m\right\}  =\left[  m\right]  $).

The $p+\left(  n-p\right)  $ elements $j_{1},j_{2},\ldots,j_{p},k_{1}%
,k_{2},\ldots,k_{n-p}$ belong to $\left[  m\right]  $ (since the $p$ elements
$j_{1},j_{2},\ldots,j_{p}$ belong to $\left[  m\right]  $, and since the $n-p$
elements $k_{1},k_{2},\ldots,k_{n-p}$ belong to $\left[  m\right]  $). In
other words, $\left(  j_{1},j_{2},\ldots,j_{p},k_{1},k_{2},\ldots
,k_{n-p}\right)  \in\left[  m\right]  ^{p+\left(  n-p\right)  }$. Thus,%
\[
\left(  j_{1},j_{2},\ldots,j_{p},k_{1},k_{2},\ldots,k_{n-p}\right)  \in\left[
m\right]  ^{p+\left(  n-p\right)  }=\left[  m\right]  ^{n}%
\]
(since $p+\left(  n-p\right)  =n$). In other words, $\left(  j_{1}%
,j_{2},\ldots,j_{p},k_{1},k_{2},\ldots,k_{n-p}\right)  $ is an $n$-tuple of
elements of $\left[  m\right]  $. Write this $n$-tuple $\left(  j_{1}%
,j_{2},\ldots,j_{p},k_{1},k_{2},\ldots,k_{n-p}\right)  $ in the form $\left(
\gamma_{1},\gamma_{2},\ldots,\gamma_{n}\right)  $. Thus, $\left(  \gamma
_{1},\gamma_{2},\ldots,\gamma_{n}\right)  =\left(  j_{1},j_{2},\ldots
,j_{p},k_{1},k_{2},\ldots,k_{n-p}\right)  $. In other words,%
\begin{equation}
\gamma_{i}=%
\begin{cases}
j_{i}, & \text{if }i\leq p;\\
k_{i-p}, & \text{if }i>p
\end{cases}
\ \ \ \ \ \ \ \ \ \ \text{for every }i\in\left\{  1,2,\ldots,n\right\}  .
\label{sol.det.laplace-multi.0r.a.def-kappa}%
\end{equation}
Moreover, $\left(  \gamma_{1},\gamma_{2},\ldots,\gamma_{n}\right)  \in\left[
m\right]  ^{n}$ (since $\left(  \gamma_{1},\gamma_{2},\ldots,\gamma
_{n}\right)  $ is an $n$-tuple of elements of $\left[  m\right]  $). Hence,
$\gamma_{i}\in\left[  m\right]  $ for each $i\in\left\{  1,2,\ldots,n\right\}
$. In other words, $\gamma_{i}\in\left[  m\right]  $ for each $i\in\left[
n\right]  $ (since $\left[  n\right]  =\left\{  1,2,\ldots,n\right\}  $).
Thus, we can define a map $\gamma:\left[  n\right]  \rightarrow\left[
m\right]  $ by%
\[
\left(  \gamma\left(  i\right)  =\gamma_{i}\text{ for every }i\in\left[
n\right]  \right)  .
\]
Consider this $\gamma$.
\end{verlong}

We have%
\begin{equation}
\gamma\left(  i\right)  =j_{i}\ \ \ \ \ \ \ \ \ \ \text{for every }%
i\in\left\{  1,2,\ldots,p\right\}  . \label{sol.det.laplace-multi.0r.a.kappa1}%
\end{equation}

\begin{vershort}
[\textit{Proof of (\ref{sol.det.laplace-multi.0r.a.kappa1}):} Let
$i\in\left\{  1,2,\ldots,p\right\}  $. Then, $i\in\left\{  1,2,\ldots
,p\right\}  \subseteq\left\{  1,2,\ldots,n\right\}  =\left[  n\right]  $.
Hence, the definition of $\gamma$ yields $\gamma\left(  i\right)  =\gamma
_{i}=j_{i}$ (by (\ref{sol.det.laplace-multi.0r.short.2a})). This proves
(\ref{sol.det.laplace-multi.0r.a.kappa1}).]
\end{vershort}

\begin{verlong}
[\textit{Proof of (\ref{sol.det.laplace-multi.0r.a.kappa1}):} Let
$i\in\left\{  1,2,\ldots,p\right\}  $. Then, $i\in\left\{  1,2,\ldots
,p\right\}  \subseteq\left\{  1,2,\ldots,n\right\}  =\left[  n\right]  $.
Hence, $\gamma\left(  i\right)  =\gamma_{i}$ (by the definition of $\gamma$).
Thus,%
\begin{align*}
\gamma\left(  i\right)   &  =\gamma_{i}=%
\begin{cases}
j_{i}, & \text{if }i\leq p;\\
k_{i-p}, & \text{if }i>p
\end{cases}
\ \ \ \ \ \ \ \ \ \ \left(  \text{by the definition of }\gamma_{i}\right) \\
&  =j_{i}\ \ \ \ \ \ \ \ \ \ \left(  \text{since }i\leq p\text{ (since }%
i\in\left\{  1,2,\ldots,p\right\}  \text{)}\right)  .
\end{align*}
This proves (\ref{sol.det.laplace-multi.0r.a.kappa1}).]
\end{verlong}

Furthermore,%
\begin{equation}
\gamma\left(  p+i\right)  =k_{i}\ \ \ \ \ \ \ \ \ \ \text{for every }%
i\in\left\{  1,2,\ldots,n-p\right\}  .
\label{sol.det.laplace-multi.0r.a.kappa2}%
\end{equation}

\begin{vershort}
[\textit{Proof of (\ref{sol.det.laplace-multi.0r.a.kappa2}):} Let
$i\in\left\{  1,2,\ldots,n-p\right\}  $. Thus, $p+i\in\left\{  p+1,p+2,\ldots
,n\right\}  \subseteq\left\{  1,2,\ldots,n\right\}  =\left[  n\right]  $.
Hence, the definition of $\gamma$ yields $\gamma\left(  p+i\right)
=\gamma_{p+i}=k_{i}$ (by (\ref{sol.det.laplace-multi.0r.short.2b})). This
proves (\ref{sol.det.laplace-multi.0r.a.kappa2}).]
\end{vershort}

\begin{verlong}
[\textit{Proof of (\ref{sol.det.laplace-multi.0r.a.kappa2}):} Let
$i\in\left\{  1,2,\ldots,n-p\right\}  $. Thus, $i\geq1>0$. But $i\in\left\{
1,2,\ldots,n-p\right\}  $, and thus%
\begin{align*}
p+i  &  \in\left\{  p+1,p+2,\ldots,p+\left(  n-p\right)  \right\}  =\left\{
p+1,p+2,\ldots,n\right\} \\
&  \ \ \ \ \ \ \ \ \ \ \left(  \text{since }p+\left(  n-p\right)  =n\right) \\
&  \subseteq\left\{  1,2,\ldots,n\right\}  \ \ \ \ \ \ \ \ \ \ \left(
\text{since }\underbrace{p}_{\geq0}+1\geq1\right) \\
&  =\left[  n\right]  .
\end{align*}
Hence, $\gamma\left(  p+i\right)  =\gamma_{p+i}$ (by the definition of
$\gamma$). Hence,%
\begin{align*}
\gamma\left(  p+i\right)   &  =\gamma_{p+i}=%
\begin{cases}
j_{p+i}, & \text{if }p+i\leq p;\\
k_{\left(  p+i\right)  -p}, & \text{if }p+i>p
\end{cases}
\ \ \ \ \ \ \ \ \ \ \left(  \text{by the definition of }\gamma_{p+i}\right) \\
&  =k_{\left(  p+i\right)  -p}\ \ \ \ \ \ \ \ \ \ \left(  \text{since
}p+\underbrace{i}_{>0}>p\right) \\
&  =k_{i}\ \ \ \ \ \ \ \ \ \ \left(  \text{since }\left(  p+i\right)
-p=i\right)  .
\end{align*}
This proves (\ref{sol.det.laplace-multi.0r.a.kappa2}).]
\end{verlong}

\begin{vershort}
We have%
\begin{equation}
w\left(  J\right)  =\left(  \gamma\left(  1\right)  ,\gamma\left(  2\right)
,\ldots,\gamma\left(  p\right)  \right)
\label{sol.det.laplace-multi.0r.a.short.wR=}%
\end{equation}
\footnote{\textit{Proof of (\ref{sol.det.laplace-multi.0r.a.short.wR=}):} The
lists $w\left(  J\right)  $ and $\left(  j_{1},j_{2},\ldots,j_{p}\right)  $
must be identical (since each of them is the list of all elements of $J$ in
increasing order (with no repetitions)). In other words, $w\left(  J\right)
=\left(  j_{1},j_{2},\ldots,j_{p}\right)  $. But
(\ref{sol.det.laplace-multi.0r.a.kappa1}) yields $\left(  \gamma\left(
1\right)  ,\gamma\left(  2\right)  ,\ldots,\gamma\left(  p\right)  \right)
=\left(  j_{1},j_{2},\ldots,j_{p}\right)  $. Comparing this with $w\left(
J\right)  =\left(  j_{1},j_{2},\ldots,j_{p}\right)  $, we obtain $w\left(
J\right)  =\left(  \gamma\left(  1\right)  ,\gamma\left(  2\right)
,\ldots,\gamma\left(  p\right)  \right)  $.} and%
\begin{equation}
w\left(  K\right)  =\left(  \gamma\left(  p+1\right)  ,\gamma\left(
p+2\right)  ,\ldots,\gamma\left(  n\right)  \right)
\label{sol.det.laplace-multi.0r.a.short.wPs=}%
\end{equation}
\footnote{\textit{Proof of (\ref{sol.det.laplace-multi.0r.a.short.wPs=}):} The
lists $w\left(  K\right)  $ and $\left(  k_{1},k_{2},\ldots,k_{n-p}\right)  $
must be identical (since each of them is the list of all elements of $K$ in
increasing order (with no repetitions)). In other words, $w\left(  K\right)
=\left(  k_{1},k_{2},\ldots,k_{n-p}\right)  $. But
(\ref{sol.det.laplace-multi.0r.a.kappa1}) yields $\left(  \gamma\left(
p+1\right)  ,\gamma\left(  p+2\right)  ,\ldots,\gamma\left(  n\right)
\right)  =\left(  k_{1},k_{2},\ldots,k_{n-p}\right)  $. Comparing this with
$w\left(  K\right)  =\left(  k_{1},k_{2},\ldots,k_{n-p}\right)  $, we obtain
$w\left(  K\right)  =\left(  \gamma\left(  p+1\right)  ,\gamma\left(
p+2\right)  ,\ldots,\gamma\left(  n\right)  \right)  $.}.
\end{vershort}

\begin{verlong}
We have%
\begin{equation}
w\left(  J\right)  =\left(  \gamma\left(  1\right)  ,\gamma\left(  2\right)
,\ldots,\gamma\left(  p\right)  \right)
\label{sol.det.laplace-multi.0r.a.wR=}%
\end{equation}
\footnote{\textit{Proof of (\ref{sol.det.laplace-multi.0r.a.wR=}):} Recall
that $w\left(  J\right)  $ is the list of all elements of $J$ in increasing
order (with no repetitions) (by the definition of $w\left(  J\right)  $).
Thus,%
\begin{align}
w\left(  J\right)   &  =\left(  \text{the list of all elements of }J\text{ in
increasing order (with no repetitions)}\right) \nonumber\\
&  =\left(  j_{1},j_{2},\ldots,j_{p}\right)
\label{sol.det.laplace-multi.0r.a.wR=.pf.1}%
\end{align}
(since $\left(  j_{1},j_{2},\ldots,j_{p}\right)  $ is the list of all elements
of $J$ in increasing order (with no repetitions)). Now,
(\ref{sol.det.laplace-multi.0r.a.kappa1}) shows that $\gamma\left(  i\right)
=j_{i}$ for every $i\in\left\{  1,2,\ldots,p\right\}  $. In other words,
$\left(  \gamma\left(  1\right)  ,\gamma\left(  2\right)  ,\ldots
,\gamma\left(  p\right)  \right)  =\left(  j_{1},j_{2},\ldots,j_{p}\right)  $.
Comparing this with (\ref{sol.det.laplace-multi.0r.a.wR=.pf.1}), we obtain
$w\left(  J\right)  =\left(  \gamma\left(  1\right)  ,\gamma\left(  2\right)
,\ldots,\gamma\left(  p\right)  \right)  $. This proves
(\ref{sol.det.laplace-multi.0r.a.wR=}).} and%
\begin{equation}
w\left(  K\right)  =\left(  \gamma\left(  p+1\right)  ,\gamma\left(
p+2\right)  ,\ldots,\gamma\left(  n\right)  \right)
\label{sol.det.laplace-multi.0r.a.wPs=}%
\end{equation}
\footnote{\textit{Proof of (\ref{sol.det.laplace-multi.0r.a.wPs=}):} Recall
that $w\left(  K\right)  $ is the list of all elements of $K$ in increasing
order (with no repetitions) (by the definition of $w\left(  K\right)  $).
Thus,%
\begin{align}
w\left(  K\right)   &  =\left(  \text{the list of all elements of }K\text{ in
increasing order (with no repetitions)}\right) \nonumber\\
&  =\left(  k_{1},k_{2},\ldots,k_{n-p}\right)
\label{sol.det.laplace-multi.0r.a.wPs=.pf.1}%
\end{align}
(since $\left(  k_{1},k_{2},\ldots,k_{n-p}\right)  $ is the list of all
elements of $K$ in increasing order (with no repetitions)). Now,
(\ref{sol.det.laplace-multi.0r.a.kappa2}) shows that $\gamma\left(
p+i\right)  =k_{i}$ for every $i\in\left\{  1,2,\ldots,n-p\right\}  $. In
other words, $\left(  \gamma\left(  p+1\right)  ,\gamma\left(  p+2\right)
,\ldots,\gamma\left(  p+\left(  n-p\right)  \right)  \right)  =\left(
k_{1},k_{2},\ldots,k_{n-p}\right)  $. Comparing this with
(\ref{sol.det.laplace-multi.0r.a.wPs=.pf.1}), we obtain
\begin{align*}
w\left(  K\right)   &  =\left(  \gamma\left(  p+1\right)  ,\gamma\left(
p+2\right)  ,\ldots,\gamma\left(  p+\left(  n-p\right)  \right)  \right) \\
&  =\left(  \gamma\left(  p+1\right)  ,\gamma\left(  p+2\right)
,\ldots,\gamma\left(  n\right)  \right)  \ \ \ \ \ \ \ \ \ \ \left(
\text{since }p+\left(  n-p\right)  =n\right)  .
\end{align*}
This proves (\ref{sol.det.laplace-multi.0r.a.wPs=}).}.
\end{verlong}

\begin{vershort}
Furthermore, there exist two distinct elements $u$ and $v$ of $\left[
n\right]  $ satisfying $\gamma\left(  u\right)  =\gamma\left(  v\right)
$\ \ \ \ \footnote{\textit{Proof.} We have $J\cap K\neq\varnothing$. In other
words, there exists some $h\in J\cap K$. Consider this $h$.
\par
Recall that $\left(  j_{1},j_{2},\ldots,j_{p}\right)  $ is the list of all
elements of $J$ in increasing order (with no repetitions). Thus, $J=\left\{
j_{1},j_{2},\ldots,j_{p}\right\}  $. Similarly, $K=\left\{  k_{1},k_{2}%
,\ldots,k_{n-p}\right\}  $.
\par
We have $h\in J\cap K\subseteq J=\left\{  j_{1},j_{2},\ldots,j_{p}\right\}  $.
In other words, there exists some $x\in\left\{  1,2,\ldots,p\right\}  $ such
that $h=j_{x}$. Consider this $x$. We have $x\in\left\{  1,2,\ldots,p\right\}
\subseteq\left\{  1,2,\ldots,n\right\}  =\left[  n\right]  $.
\par
We have $h\in J\cap K\subseteq K=\left\{  k_{1},k_{2},\ldots,k_{n-p}\right\}
$. In other words, there exists some $y\in\left\{  1,2,\ldots,n-p\right\}  $
such that $h=k_{y}$. Consider this $y$. We have $y\in\left\{  1,2,\ldots
,n-p\right\}  $ and thus $p+y\in\left\{  p+1,p+2,\ldots,n\right\}
\subseteq\left\{  1,2,\ldots,n\right\}  =\left[  n\right]  $.
\par
From $x\in\left\{  1,2,\ldots,p\right\}  $, we obtain $p\geq x$. From
$y\in\left\{  1,2,\ldots,n-p\right\}  $, we obtain $y>0$, so that
$p+\underbrace{y}_{>0}>p\geq x$. Hence, $p+y\neq x$. Thus, the numbers $x$ and
$p+y$ are distinct.
\par
Both $x$ and $p+y$ are elements of $\left[  n\right]  $ (since $x\in\left[
n\right]  $ and $p+y\in\left[  n\right]  $). Hence, $x$ and $p+y$ are two
distinct elements of $\left[  n\right]  $.
\par
Applying (\ref{sol.det.laplace-multi.0r.a.kappa1}) to $i=x$, we obtain
$\gamma\left(  x\right)  =j_{x}=h$ (since $h=j_{x}$).
\par
Applying (\ref{sol.det.laplace-multi.0r.a.kappa2}) to $i=y$, we obtain
$\gamma\left(  p+y\right)  =k_{y}=h$ (since $h=k_{y}$).
\par
Thus, $\gamma\left(  x\right)  =h=\gamma\left(  p+y\right)  $. Hence, there
exist two distinct elements $u$ and $v$ of $\left[  n\right]  $ satisfying
$\gamma\left(  u\right)  =\gamma\left(  v\right)  $ (namely, $u=x$ and
$v=p+y$).}.
\end{vershort}

\begin{verlong}
Furthermore, there exist two distinct elements $u$ and $v$ of $\left[
n\right]  $ satisfying $\gamma\left(  u\right)  =\gamma\left(  v\right)
$\ \ \ \ \footnote{\textit{Proof.} We have $J\cap K\neq\varnothing$. Thus, the
set $J\cap K$ is nonempty, and therefore contains at least one element. In
other words, there exists some $h\in J\cap K$. Consider this $h$.
\par
Recall that $\left(  j_{1},j_{2},\ldots,j_{p}\right)  $ is the list of all
elements of $J$ in increasing order (with no repetitions). Thus, $J=\left\{
j_{1},j_{2},\ldots,j_{p}\right\}  $.
\par
Recall that $\left(  k_{1},k_{2},\ldots,k_{n-p}\right)  $ is the list of all
elements of $K$ in increasing order (with no repetitions). Thus, $K=\left\{
k_{1},k_{2},\ldots,k_{n-p}\right\}  $.
\par
We have $h\in J\cap K\subseteq J=\left\{  j_{1},j_{2},\ldots,j_{p}\right\}  $.
In other words, there exists some $x\in\left\{  1,2,\ldots,p\right\}  $ such
that $h=j_{x}$. Consider this $x$. We have $x\in\left\{  1,2,\ldots,p\right\}
\subseteq\left\{  1,2,\ldots,n\right\}  $.
\par
We have $h\in J\cap K\subseteq K=\left\{  k_{1},k_{2},\ldots,k_{n-p}\right\}
$. In other words, there exists some $y\in\left\{  1,2,\ldots,n-p\right\}  $
such that $h=k_{y}$. Consider this $y$. We have $y\in\left\{  1,2,\ldots
,n-p\right\}  $ and thus
\begin{align*}
p+y  &  \in\left\{  p+1,p+2,\ldots,p+\left(  n-p\right)  \right\}  =\left\{
p+1,p+2,\ldots,n\right\}  \ \ \ \ \ \ \ \ \ \ \left(  \text{since }p+\left(
n-p\right)  =n\right) \\
&  \subseteq\left\{  1,2,\ldots,n\right\}  .
\end{align*}
\par
From $x\in\left\{  1,2,\ldots,p\right\}  $, we obtain $x\leq p$, and thus
$p\geq x$. From $y\in\left\{  1,2,\ldots,n-p\right\}  $, we obtain $y\geq1>0$,
so that $p+\underbrace{y}_{>0}>p\geq x$. Hence, $p+y\neq x$. Thus, the numbers
$x$ and $p+y$ are distinct.
\par
Both $x$ and $p+y$ are elements of $\left\{  1,2,\ldots,n\right\}  $ (since
$x\in\left\{  1,2,\ldots,n\right\}  $ and $p+y\in\left\{  1,2,\ldots
,n\right\}  $), and thus are elements of $\left[  n\right]  $ (since $\left\{
1,2,\ldots,n\right\}  =\left[  n\right]  $). Hence, $x$ and $p+y$ are two
distinct elements of $\left[  n\right]  $.
\par
Applying (\ref{sol.det.laplace-multi.0r.a.kappa1}) to $i=x$, we obtain
$\gamma\left(  x\right)  =j_{x}=h$ (since $h=j_{x}$).
\par
Applying (\ref{sol.det.laplace-multi.0r.a.kappa2}) to $i=y$, we obtain
$\gamma\left(  p+y\right)  =k_{y}=h$ (since $h=k_{y}$).
\par
Thus, $\gamma\left(  x\right)  =h=\gamma\left(  p+y\right)  $. Hence, there
exist two distinct elements $u$ and $v$ of $\left[  n\right]  $ satisfying
$\gamma\left(  u\right)  =\gamma\left(  v\right)  $ (namely, $u=x$ and
$v=p+y$). Qed.}.
\end{verlong}

\begin{vershort}
Write the $m\times n$-matrix $A$ in the form $A=\left(  a_{i,j}\right)
_{1\leq i\leq m,\ 1\leq j\leq n}$. Define an $n\times n$-matrix $A_{\gamma}$
by $A_{\gamma}=\left(  a_{\gamma\left(  i\right)  ,j}\right)  _{1\leq i\leq
n,\ 1\leq j\leq n}$. Then, the matrix $A_{\gamma}$ has two equal
rows\footnote{\textit{Proof.} We have shown that there exist two distinct
elements $u$ and $v$ of $\left[  n\right]  $ satisfying $\gamma\left(
u\right)  =\gamma\left(  v\right)  $. Consider these $u$ and $v$.
\par
We have $A_{\gamma}=\left(  a_{\gamma\left(  i\right)  ,j}\right)  _{1\leq
i\leq n,\ 1\leq j\leq n}$. Thus,%
\[
\left(  \text{the }u\text{-th row of the matrix }A_{\gamma}\right)  =\left(
a_{\gamma\left(  u\right)  ,1},a_{\gamma\left(  u\right)  ,2},\ldots
,a_{\gamma\left(  u\right)  ,n}\right)  =\left(  a_{\gamma\left(  v\right)
,1},a_{\gamma\left(  v\right)  ,2},\ldots,a_{\gamma\left(  v\right)
,n}\right)
\]
(since $\gamma\left(  u\right)  =\gamma\left(  v\right)  $). Comparing this
with%
\begin{align*}
&  \left(  \text{the }v\text{-th row of the matrix }A_{\gamma}\right) \\
&  =\left(  a_{\gamma\left(  v\right)  ,1},a_{\gamma\left(  v\right)
,2},\ldots,a_{\gamma\left(  v\right)  ,n}\right)  \ \ \ \ \ \ \ \ \ \ \left(
\text{since }A_{\gamma}=\left(  a_{\gamma\left(  i\right)  ,j}\right)  _{1\leq
i\leq n,\ 1\leq j\leq n}\right)  ,
\end{align*}
we obtain $\left(  \text{the }u\text{-th row of the matrix }A_{\gamma}\right)
=\left(  \text{the }v\text{-th row of the matrix }A_{\gamma}\right)  $. Since
$u$ and $v$ are distinct, this shows that the matrix $A_{\gamma}$ has two
equal rows.}.
\end{vershort}

\begin{verlong}
Write the $m\times n$-matrix $A$ in the form $A=\left(  a_{i,j}\right)
_{1\leq i\leq m,\ 1\leq j\leq n}$. Define an $n\times n$-matrix $A_{\gamma}$
by $A_{\gamma}=\left(  a_{\gamma\left(  i\right)  ,j}\right)  _{1\leq i\leq
n,\ 1\leq j\leq n}$. Then, the matrix $A_{\gamma}$ has two equal
rows\footnote{\textit{Proof.} We have shown that there exist two distinct
elements $u$ and $v$ of $\left[  n\right]  $ satisfying $\gamma\left(
u\right)  =\gamma\left(  v\right)  $. Consider these $u$ and $v$.
\par
We have $u\in\left[  n\right]  =\left\{  1,2,\ldots,n\right\}  $ and
$v\in\left[  n\right]  =\left\{  1,2,\ldots,n\right\}  $. Thus, the $n\times
n$-matrix $A_{\gamma}$ has a $u$-th row and a $v$-th row.
\par
We have $A_{\gamma}=\left(  a_{\gamma\left(  i\right)  ,j}\right)  _{1\leq
i\leq n,\ 1\leq j\leq n}$. Thus,%
\begin{align*}
&  \left(  \text{the }u\text{-th row of the matrix }A_{\gamma}\right) \\
&  =\left(  \text{the }u\text{-th row of the matrix }\left(  a_{\gamma\left(
i\right)  ,j}\right)  _{1\leq i\leq n,\ 1\leq j\leq n}\right) \\
&  =\left(  a_{\gamma\left(  u\right)  ,1},a_{\gamma\left(  u\right)
,2},\ldots,a_{\gamma\left(  u\right)  ,n}\right)  =\left(  a_{\gamma\left(
v\right)  ,1},a_{\gamma\left(  v\right)  ,2},\ldots,a_{\gamma\left(  v\right)
,n}\right)
\end{align*}
(since $\gamma\left(  u\right)  =\gamma\left(  v\right)  $). Comparing this
with%
\begin{align*}
&  \left(  \text{the }v\text{-th row of the matrix }\underbrace{A_{\gamma}%
}_{=\left(  a_{\gamma\left(  i\right)  ,j}\right)  _{1\leq i\leq n,\ 1\leq
j\leq n}}\right) \\
&  =\left(  \text{the }v\text{-th row of the matrix }\left(  a_{\gamma\left(
i\right)  ,j}\right)  _{1\leq i\leq n,\ 1\leq j\leq n}\right) \\
&  =\left(  a_{\gamma\left(  v\right)  ,1},a_{\gamma\left(  v\right)
,2},\ldots,a_{\gamma\left(  v\right)  ,n}\right)  ,
\end{align*}
we obtain $\left(  \text{the }u\text{-th row of the matrix }A_{\gamma}\right)
=\left(  \text{the }v\text{-th row of the matrix }A_{\gamma}\right)  $. In
other words, the $u$-th row and the $v$-th row of the matrix $A_{\gamma}$ are
equal. Since $u$ and $v$ are distinct, this shows that the matrix $A_{\gamma}$
has two equal rows. Qed.}.
\end{verlong}

Hence, Exercise \ref{exe.ps4.6} \textbf{(e)} (applied to $A_{\gamma}$ instead
of $A$) shows that
\begin{equation}
\det\left(  A_{\gamma}\right)  =0. \label{sol.det.laplace-multi.0r.a.detAk}%
\end{equation}

Now, every subset $Q$ of $\left\{  1,2,\ldots,n\right\}  $ satisfies%
\begin{equation}
\operatorname*{sub}\nolimits_{w\left(  J\right)  }^{w\left(  Q\right)
}A=\operatorname*{sub}\nolimits_{\left(  1,2,\ldots,p\right)  }^{w\left(
Q\right)  }\left(  A_{\gamma}\right)  \label{sol.det.laplace-multi.0r.a.sub1}%
\end{equation}
\footnote{\textit{Proof of (\ref{sol.det.laplace-multi.0r.a.sub1}):} Let $Q$
be a subset of $\left\{  1,2,\ldots,n\right\}  $. Write the list $w\left(
Q\right)  $ in the form $w\left(  Q\right)  =\left(  q_{1},q_{2}%
,\ldots,q_{\ell}\right)  $ for some $\ell\in\mathbb{N}$.
\par
From $w\left(  J\right)  =\left(  \gamma\left(  1\right)  ,\gamma\left(
2\right)  ,\ldots,\gamma\left(  p\right)  \right)  $ and $w\left(  Q\right)
=\left(  q_{1},q_{2},\ldots,q_{\ell}\right)  $, we obtain%
\begin{align}
\operatorname*{sub}\nolimits_{w\left(  J\right)  }^{w\left(  Q\right)  }A  &
=\operatorname*{sub}\nolimits_{\left(  \gamma\left(  1\right)  ,\gamma\left(
2\right)  ,\ldots,\gamma\left(  p\right)  \right)  }^{\left(  q_{1}%
,q_{2},\ldots,q_{\ell}\right)  }A=\operatorname*{sub}\nolimits_{\gamma\left(
1\right)  ,\gamma\left(  2\right)  ,\ldots,\gamma\left(  p\right)  }%
^{q_{1},q_{2},\ldots,q_{\ell}}A=\left(  a_{\gamma\left(  x\right)  ,q_{y}%
}\right)  _{1\leq x\leq p,\ 1\leq y\leq\ell}%
\label{sol.det.laplace-multi.0r.a.sub1.pf.1}\\
&  \ \ \ \ \ \ \ \ \ \ \left(  \text{by the definition of }\operatorname*{sub}%
\nolimits_{\gamma\left(  1\right)  ,\gamma\left(  2\right)  ,\ldots
,\gamma\left(  p\right)  }^{q_{1},q_{2},\ldots,q_{\ell}}A\text{, since
}A=\left(  a_{i,j}\right)  _{1\leq i\leq m,\ 1\leq j\leq n}\right)  .\nonumber
\end{align}
\par
On the other hand, from $w\left(  Q\right)  =\left(  q_{1},q_{2}%
,\ldots,q_{\ell}\right)  $, we obtain%
\[
\operatorname*{sub}\nolimits_{\left(  1,2,\ldots,p\right)  }^{w\left(
Q\right)  }\left(  A_{\gamma}\right)  =\operatorname*{sub}\nolimits_{\left(
1,2,\ldots,p\right)  }^{\left(  q_{1},q_{2},\ldots,q_{\ell}\right)  }\left(
A_{\gamma}\right)  =\operatorname*{sub}\nolimits_{1,2,\ldots,p}^{q_{1}%
,q_{2},\ldots,q_{\ell}}\left(  A_{\gamma}\right)  =\left(  a_{\gamma\left(
x\right)  ,q_{y}}\right)  _{1\leq x\leq p,\ 1\leq y\leq\ell}%
\]
(by the definition of $\operatorname*{sub}\nolimits_{1,2,\ldots,p}%
^{q_{1},q_{2},\ldots,q_{\ell}}\left(  A_{\gamma}\right)  $, since $A_{\gamma
}=\left(  a_{\gamma\left(  i\right)  ,j}\right)  _{1\leq i\leq n,\ 1\leq j\leq
n}$). Comparing this with (\ref{sol.det.laplace-multi.0r.a.sub1.pf.1}), we
obtain $\operatorname*{sub}\nolimits_{w\left(  J\right)  }^{w\left(  Q\right)
}A=\operatorname*{sub}\nolimits_{\left(  1,2,\ldots,p\right)  }^{w\left(
Q\right)  }\left(  A_{\gamma}\right)  $. This proves
(\ref{sol.det.laplace-multi.0r.a.sub1}).} and%
\begin{equation}
\operatorname*{sub}\nolimits_{w\left(  K\right)  }^{w\left(  \widetilde{Q}%
\right)  }A=\operatorname*{sub}\nolimits_{\left(  p+1,p+2,\ldots,n\right)
}^{w\left(  \widetilde{Q}\right)  }\left(  A_{\gamma}\right)
\label{sol.det.laplace-multi.0r.a.sub2}%
\end{equation}
\footnote{\textit{Proof of (\ref{sol.det.laplace-multi.0r.a.sub2}):} Let $Q$
be a subset of $\left\{  1,2,\ldots,n\right\}  $. Write the list $w\left(
\widetilde{Q}\right)  $ in the form $w\left(  \widetilde{Q}\right)  =\left(
q_{1},q_{2},\ldots,q_{\ell}\right)  $ for some $\ell\in\mathbb{N}$.
\par
From $w\left(  K\right)  =\left(  \gamma\left(  p+1\right)  ,\gamma\left(
p+2\right)  ,\ldots,\gamma\left(  n\right)  \right)  $ and $w\left(
\widetilde{Q}\right)  =\left(  q_{1},q_{2},\ldots,q_{\ell}\right)  $, we
obtain%
\begin{align}
\operatorname*{sub}\nolimits_{w\left(  K\right)  }^{w\left(  \widetilde{Q}%
\right)  }A  &  =\operatorname*{sub}\nolimits_{\left(  \gamma\left(
p+1\right)  ,\gamma\left(  p+2\right)  ,\ldots,\gamma\left(  n\right)
\right)  }^{\left(  q_{1},q_{2},\ldots,q_{\ell}\right)  }A=\operatorname*{sub}%
\nolimits_{\gamma\left(  p+1\right)  ,\gamma\left(  p+2\right)  ,\ldots
,\gamma\left(  n\right)  }^{q_{1},q_{2},\ldots,q_{\ell}}A=\left(
a_{\gamma\left(  p+x\right)  ,q_{y}}\right)  _{1\leq x\leq n-p,\ 1\leq
y\leq\ell}\label{sol.det.laplace-multi.0r.a.sub2.pf.1}\\
&  \ \ \ \ \ \ \ \ \ \ \left(  \text{by the definition of }\operatorname*{sub}%
\nolimits_{\gamma\left(  p+1\right)  ,\gamma\left(  p+2\right)  ,\ldots
,\gamma\left(  n\right)  }^{q_{1},q_{2},\ldots,q_{\ell}}A\text{, since
}A=\left(  a_{i,j}\right)  _{1\leq i\leq m,\ 1\leq j\leq n}\right)  .\nonumber
\end{align}
\par
On the other hand, from $w\left(  \widetilde{Q}\right)  =\left(  q_{1}%
,q_{2},\ldots,q_{\ell}\right)  $, we obtain%
\[
\operatorname*{sub}\nolimits_{\left(  p+1,p+2,\ldots,n\right)  }^{w\left(
\widetilde{Q}\right)  }\left(  A_{\gamma}\right)  =\operatorname*{sub}%
\nolimits_{\left(  p+1,p+2,\ldots,n\right)  }^{\left(  q_{1},q_{2}%
,\ldots,q_{\ell}\right)  }\left(  A_{\gamma}\right)  =\operatorname*{sub}%
\nolimits_{p+1,p+2,\ldots,n}^{q_{1},q_{2},\ldots,q_{\ell}}\left(  A_{\gamma
}\right)  =\left(  a_{\gamma\left(  p+x\right)  ,q_{y}}\right)  _{1\leq x\leq
n-p,\ 1\leq y\leq\ell}%
\]
(by the definition of $\operatorname*{sub}\nolimits_{p+1,p+2,\ldots,n}%
^{q_{1},q_{2},\ldots,q_{\ell}}\left(  A_{\gamma}\right)  $, since $A_{\gamma
}=\left(  a_{\gamma\left(  i\right)  ,j}\right)  _{1\leq i\leq n,\ 1\leq j\leq
n}$). Comparing this with (\ref{sol.det.laplace-multi.0r.a.sub2.pf.1}), we
obtain $\operatorname*{sub}\nolimits_{w\left(  K\right)  }^{w\left(
\widetilde{Q}\right)  }A=\operatorname*{sub}\nolimits_{\left(  p+1,p+2,\ldots
,n\right)  }^{w\left(  \widetilde{Q}\right)  }\left(  A_{\gamma}\right)  $.
This proves (\ref{sol.det.laplace-multi.0r.a.sub2}).}.

But $\left\{  1,2,\ldots,p\right\}  \subseteq\left\{  1,2,\ldots,n\right\}  $.
In other words, $\left\{  1,2,\ldots,p\right\}  $ is a subset of $\left\{
1,2,\ldots,n\right\}  $. Thus, we can define a subset $P^{\prime}$ of
$\left\{  1,2,\ldots,n\right\}  $ by $P^{\prime}=\left\{  1,2,\ldots
,p\right\}  $. Consider this $P^{\prime}$. From $P^{\prime}=\left\{
1,2,\ldots,p\right\}  $, we obtain $\left\vert P^{\prime}\right\vert
=\left\vert \left\{  1,2,\ldots,p\right\}  \right\vert =p=\left\vert
J\right\vert $.

We have $w\left(  P^{\prime}\right)  =\left(  1,2,\ldots,p\right)
$\ \ \ \ \footnote{\textit{Proof.} Recall that $w\left(  P^{\prime}\right)  $
is the list of all elements of $P^{\prime}$ in increasing order (with no
repetitions) (by the definition of $w\left(  P^{\prime}\right)  $). Thus,%
\begin{align*}
w\left(  P^{\prime}\right)   &  =\left(  \text{the list of all elements of
}\underbrace{P^{\prime}}_{=\left\{  1,2,\ldots,p\right\}  }\text{ in
increasing order (with no repetitions)}\right) \\
&  =\left(  \text{the list of all elements of }\left\{  1,2,\ldots,p\right\}
\text{ in increasing order (with no repetitions)}\right) \\
&  =\left(  1,2,\ldots,p\right)  .
\end{align*}
Qed.} and $w\left(  \widetilde{P^{\prime}}\right)  =\left(  p+1,p+2,\ldots
,n\right)  $\ \ \ \ \footnote{\textit{Proof.} We have%
\begin{align*}
\widetilde{P^{\prime}}  &  =\left\{  1,2,\ldots,n\right\}  \setminus
\underbrace{P^{\prime}}_{=\left\{  1,2,\ldots,p\right\}  }%
\ \ \ \ \ \ \ \ \ \ \left(  \text{by the definition of }\widetilde{P^{\prime}%
}\right) \\
&  =\left\{  1,2,\ldots,n\right\}  \setminus\left\{  1,2,\ldots,p\right\}
=\left\{  p+1,p+2,\ldots,n\right\}  \ \ \ \ \ \ \ \ \ \ \left(  \text{since
}p\in\left\{  0,1,\ldots,n\right\}  \right)  .
\end{align*}
\par
Recall that $w\left(  \widetilde{P^{\prime}}\right)  $ is the list of all
elements of $\widetilde{P^{\prime}}$ in increasing order (with no repetitions)
(by the definition of $w\left(  \widetilde{P^{\prime}}\right)  $). Thus,%
\begin{align*}
&  w\left(  \widetilde{P^{\prime}}\right) \\
&  =\left(  \text{the list of all elements of }%
\underbrace{\widetilde{P^{\prime}}}_{=\left\{  p+1,p+2,\ldots,n\right\}
}\text{ in increasing order (with no repetitions)}\right) \\
&  =\left(  \text{the list of all elements of }\left\{  p+1,p+2,\ldots
,n\right\}  \text{ in increasing order (with no repetitions)}\right) \\
&  =\left(  p+1,p+2,\ldots,n\right)  .
\end{align*}
Qed.}.

Now, Theorem \ref{thm.det.laplace-multi} \textbf{(a)} (applied to $A_{\gamma}$
and $P^{\prime}$ instead of $A$ and $P$) yields%
\begin{align*}
&  \det\left(  A_{\gamma}\right) \\
&  =\underbrace{\sum_{\substack{Q\subseteq\left\{  1,2,\ldots,n\right\}
;\\\left\vert Q\right\vert =\left\vert P^{\prime}\right\vert }}}%
_{\substack{=\sum_{\substack{Q\subseteq\left\{  1,2,\ldots,n\right\}
;\\\left\vert Q\right\vert =\left\vert J\right\vert }}\\\text{(since
}\left\vert P^{\prime}\right\vert =\left\vert J\right\vert \text{)}}}\left(
-1\right)  ^{\sum P^{\prime}+\sum Q}\det\left(
\underbrace{\operatorname*{sub}\nolimits_{w\left(  P^{\prime}\right)
}^{w\left(  Q\right)  }\left(  A_{\gamma}\right)  }%
_{\substack{=\operatorname*{sub}\nolimits_{\left(  1,2,\ldots,p\right)
}^{w\left(  Q\right)  }\left(  A_{\gamma}\right)  \\\text{(since }w\left(
P^{\prime}\right)  =\left(  1,2,\ldots,p\right)  \text{)}}}\right)
\det\left(  \underbrace{\operatorname*{sub}\nolimits_{w\left(
\widetilde{P^{\prime}}\right)  }^{w\left(  \widetilde{Q}\right)  }\left(
A_{\gamma}\right)  }_{\substack{=\operatorname*{sub}\nolimits_{\left(
p+1,p+2,\ldots,n\right)  }^{w\left(  \widetilde{Q}\right)  }\left(  A_{\gamma
}\right)  \\\text{(since }w\left(  \widetilde{P^{\prime}}\right)  =\left(
p+1,p+2,\ldots,n\right)  \text{)}}}\right) \\
&  =\sum_{\substack{Q\subseteq\left\{  1,2,\ldots,n\right\}  ;\\\left\vert
Q\right\vert =\left\vert J\right\vert }}\underbrace{\left(  -1\right)  ^{\sum
P^{\prime}+\sum Q}}_{=\left(  -1\right)  ^{\sum P^{\prime}}\left(  -1\right)
^{\sum Q}}\det\left(  \underbrace{\operatorname*{sub}\nolimits_{\left(
1,2,\ldots,p\right)  }^{w\left(  Q\right)  }\left(  A_{\gamma}\right)
}_{\substack{=\operatorname*{sub}\nolimits_{w\left(  J\right)  }^{w\left(
Q\right)  }A\\\text{(by (\ref{sol.det.laplace-multi.0r.a.sub1}))}}}\right)
\det\left(  \underbrace{\operatorname*{sub}\nolimits_{\left(  p+1,p+2,\ldots
,n\right)  }^{w\left(  \widetilde{Q}\right)  }\left(  A_{\gamma}\right)
}_{\substack{=\operatorname*{sub}\nolimits_{w\left(  K\right)  }^{w\left(
\widetilde{Q}\right)  }A\\\text{(by (\ref{sol.det.laplace-multi.0r.a.sub2}))}%
}}\right) \\
&  =\sum_{\substack{Q\subseteq\left\{  1,2,\ldots,n\right\}  ;\\\left\vert
Q\right\vert =\left\vert J\right\vert }}\left(  -1\right)  ^{\sum P^{\prime}%
}\left(  -1\right)  ^{\sum Q}\det\left(  \operatorname*{sub}%
\nolimits_{w\left(  J\right)  }^{w\left(  Q\right)  }A\right)  \det\left(
\operatorname*{sub}\nolimits_{w\left(  K\right)  }^{w\left(  \widetilde{Q}%
\right)  }A\right)  .
\end{align*}
Comparing this with (\ref{sol.det.laplace-multi.0r.a.detAk}), we obtain%
\[
0=\sum_{\substack{Q\subseteq\left\{  1,2,\ldots,n\right\}  ;\\\left\vert
Q\right\vert =\left\vert J\right\vert }}\left(  -1\right)  ^{\sum P^{\prime}%
}\left(  -1\right)  ^{\sum Q}\det\left(  \operatorname*{sub}%
\nolimits_{w\left(  J\right)  }^{w\left(  Q\right)  }A\right)  \det\left(
\operatorname*{sub}\nolimits_{w\left(  K\right)  }^{w\left(  \widetilde{Q}%
\right)  }A\right)  .
\]
Multiplying both sides of this equality by $\left(  -1\right)  ^{\sum
P^{\prime}}$, we obtain%
\begin{align*}
0  &  =\left(  -1\right)  ^{\sum P^{\prime}}\cdot\sum_{\substack{Q\subseteq
\left\{  1,2,\ldots,n\right\}  ;\\\left\vert Q\right\vert =\left\vert
J\right\vert }}\left(  -1\right)  ^{\sum P^{\prime}}\left(  -1\right)  ^{\sum
Q}\det\left(  \operatorname*{sub}\nolimits_{w\left(  J\right)  }^{w\left(
Q\right)  }A\right)  \det\left(  \operatorname*{sub}\nolimits_{w\left(
K\right)  }^{w\left(  \widetilde{Q}\right)  }A\right) \\
&  =\sum_{\substack{Q\subseteq\left\{  1,2,\ldots,n\right\}  ;\\\left\vert
Q\right\vert =\left\vert J\right\vert }}\underbrace{\left(  -1\right)  ^{\sum
P^{\prime}}\left(  -1\right)  ^{\sum P^{\prime}}}_{\substack{=\left(
-1\right)  ^{\sum P^{\prime}+\sum P^{\prime}}=1\\\text{(since }\sum P^{\prime
}+\sum P^{\prime}=2\sum P^{\prime}\text{ is even)}}}\left(  -1\right)  ^{\sum
Q}\det\left(  \operatorname*{sub}\nolimits_{w\left(  J\right)  }^{w\left(
Q\right)  }A\right)  \det\left(  \operatorname*{sub}\nolimits_{w\left(
K\right)  }^{w\left(  \widetilde{Q}\right)  }A\right) \\
&  =\sum_{\substack{Q\subseteq\left\{  1,2,\ldots,n\right\}  ;\\\left\vert
Q\right\vert =\left\vert J\right\vert }}\left(  -1\right)  ^{\sum Q}%
\det\left(  \operatorname*{sub}\nolimits_{w\left(  J\right)  }^{w\left(
Q\right)  }A\right)  \det\left(  \operatorname*{sub}\nolimits_{w\left(
K\right)  }^{w\left(  \widetilde{Q}\right)  }A\right)  .
\end{align*}
This solves Exercise \ref{exe.det.laplace-multi.0r} \textbf{(a)}.

\textbf{(b)} Let $A\in\mathbb{K}^{n\times m}$. Every subset $Q$ of $\left\{
1,2,\ldots,n\right\}  $ satisfying $\left\vert Q\right\vert =\left\vert
J\right\vert $ satisfies%
\begin{align}
&  \det\left(  \operatorname*{sub}\nolimits_{w\left(  J\right)  }^{w\left(
Q\right)  }\left(  A^{T}\right)  \right)  \det\left(  \operatorname*{sub}%
\nolimits_{w\left(  K\right)  }^{w\left(  \widetilde{Q}\right)  }\left(
A^{T}\right)  \right) \nonumber\\
&  =\det\left(  \operatorname*{sub}\nolimits_{w\left(  Q\right)  }^{w\left(
J\right)  }A\right)  \det\left(  \operatorname*{sub}\nolimits_{w\left(
\widetilde{Q}\right)  }^{w\left(  K\right)  }A\right)
\label{sol.det.laplace-multi.0r.b.1}%
\end{align}
\footnote{\textit{Proof of (\ref{sol.det.laplace-multi.0r.b.1}):} Let $Q$ be a
subset of $\left\{  1,2,\ldots,n\right\}  $ satisfying $\left\vert
Q\right\vert =\left\vert J\right\vert $. Thus, $\left\vert J\right\vert
=\left\vert Q\right\vert $. Hence, Corollary
\ref{cor.sol.det.laplace-multi.0.detsubAT} (applied to $U=J$ and $V=Q$)
yields
\begin{equation}
\det\left(  \operatorname*{sub}\nolimits_{w\left(  J\right)  }^{w\left(
Q\right)  }\left(  A^{T}\right)  \right)  =\det\left(  \operatorname*{sub}%
\nolimits_{w\left(  Q\right)  }^{w\left(  J\right)  }A\right)  .
\label{sol.det.laplace-multi.0r.b.1.pf.1a}%
\end{equation}
\par
On the other hand, the definition of $\widetilde{Q}$ yields $\widetilde{Q}%
=\left\{  1,2,\ldots,n\right\}  \setminus Q\subseteq\left\{  1,2,\ldots
,n\right\}  $. Hence, $\widetilde{Q}$ is a subset of $\left\{  1,2,\ldots
,n\right\}  $. Also, $\widetilde{Q}=\left\{  1,2,\ldots,n\right\}  \setminus
Q$ and thus
\begin{align*}
\left\vert \widetilde{Q}\right\vert  &  =\left\vert \left\{  1,2,\ldots
,n\right\}  \setminus Q\right\vert =\underbrace{\left\vert \left\{
1,2,\ldots,n\right\}  \right\vert }_{=n}-\underbrace{\left\vert Q\right\vert
}_{=\left\vert J\right\vert }\ \ \ \ \ \ \ \ \ \ \left(  \text{since }Q\text{
is a subset of }\left\{  1,2,\ldots,n\right\}  \right) \\
&  =n-\left\vert J\right\vert =\left\vert K\right\vert
\ \ \ \ \ \ \ \ \ \ \left(  \text{since }n=\left\vert J\right\vert +\left\vert
K\right\vert \right)  .
\end{align*}
Hence, $\left\vert K\right\vert =\left\vert \widetilde{Q}\right\vert $. Thus,
Corollary \ref{cor.sol.det.laplace-multi.0.detsubAT} (applied to $U=K$ and
$V=\widetilde{Q}$) yields
\begin{equation}
\det\left(  \operatorname*{sub}\nolimits_{w\left(  K\right)  }^{w\left(
\widetilde{Q}\right)  }\left(  A^{T}\right)  \right)  =\det\left(
\operatorname*{sub}\nolimits_{w\left(  \widetilde{Q}\right)  }^{w\left(
K\right)  }A\right)  . \label{sol.det.laplace-multi.0r.b.1.pf.2a}%
\end{equation}
Multiplying the equality (\ref{sol.det.laplace-multi.0r.b.1.pf.1a}) with the
equality (\ref{sol.det.laplace-multi.0r.b.1.pf.2a}), we obtain%
\[
\det\left(  \operatorname*{sub}\nolimits_{w\left(  J\right)  }^{w\left(
Q\right)  }\left(  A^{T}\right)  \right)  \det\left(  \operatorname*{sub}%
\nolimits_{w\left(  K\right)  }^{w\left(  \widetilde{Q}\right)  }\left(
A^{T}\right)  \right)  =\det\left(  \operatorname*{sub}\nolimits_{w\left(
Q\right)  }^{w\left(  J\right)  }A\right)  \det\left(  \operatorname*{sub}%
\nolimits_{w\left(  \widetilde{Q}\right)  }^{w\left(  K\right)  }A\right)  .
\]
Thus, (\ref{sol.det.laplace-multi.0r.b.1}) is proven.}.

We have $A\in\mathbb{K}^{n\times m}$, and thus $A^{T}\in\mathbb{K}^{m\times
n}$. Hence, Exercise \ref{exe.det.laplace-multi.0r} \textbf{(a)} (applied to
$A^{T}$ instead of $A$) yields%
\begin{align}
0  &  =\sum_{\substack{Q\subseteq\left\{  1,2,\ldots,n\right\}  ;\\\left\vert
Q\right\vert =\left\vert J\right\vert }}\left(  -1\right)  ^{\sum
Q}\underbrace{\det\left(  \operatorname*{sub}\nolimits_{w\left(  J\right)
}^{w\left(  Q\right)  }\left(  A^{T}\right)  \right)  \det\left(
\operatorname*{sub}\nolimits_{w\left(  K\right)  }^{w\left(  \widetilde{Q}%
\right)  }\left(  A^{T}\right)  \right)  }_{\substack{=\det\left(
\operatorname*{sub}\nolimits_{w\left(  Q\right)  }^{w\left(  J\right)
}A\right)  \det\left(  \operatorname*{sub}\nolimits_{w\left(  \widetilde{Q}%
\right)  }^{w\left(  K\right)  }A\right)  \\\text{(by
(\ref{sol.det.laplace-multi.0r.b.1}))}}}\nonumber\\
&  =\sum_{\substack{Q\subseteq\left\{  1,2,\ldots,n\right\}  ;\\\left\vert
Q\right\vert =\left\vert J\right\vert }}\left(  -1\right)  ^{\sum Q}%
\det\left(  \operatorname*{sub}\nolimits_{w\left(  Q\right)  }^{w\left(
J\right)  }A\right)  \det\left(  \operatorname*{sub}\nolimits_{w\left(
\widetilde{Q}\right)  }^{w\left(  K\right)  }A\right) \nonumber\\
&  =\sum_{\substack{P\subseteq\left\{  1,2,\ldots,n\right\}  ;\\\left\vert
P\right\vert =\left\vert J\right\vert }}\left(  -1\right)  ^{\sum P}%
\det\left(  \operatorname*{sub}\nolimits_{w\left(  P\right)  }^{w\left(
J\right)  }A\right)  \det\left(  \operatorname*{sub}\nolimits_{w\left(
\widetilde{P}\right)  }^{w\left(  K\right)  }A\right)
\label{sol.det.laplace-multi.0r.b.5}%
\end{align}
(here, we have renamed the summation index $Q$ as $P$). This solves Exercise
\ref{exe.det.laplace-multi.0r} \textbf{(b)}.
\end{proof}

\subsection{Solution to Exercise \ref{exe.det.blocktria-twisted}}

Before we solve Exercise \ref{exe.det.blocktria-twisted} using Theorem
\ref{thm.det.laplace-multi}, let us outline a quick solution for Exercise
\ref{exe.det.blocktria-twisted} \textbf{(a)} using just the definition of the
determinant (generalizing our first solution to Exercise \ref{exe.ps4.5}
\textbf{(a)}):

\begin{proof}
[First solution to Exercise \ref{exe.det.blocktria-twisted} \textbf{(a)}
(sketched).]Assume that $\left\vert P\right\vert +\left\vert Q\right\vert >n$.

Let $\sigma\in S_{n}$. We shall write $\left[  n\right]  $ for the set
$\left\{  1,2,\ldots,n\right\}  $.

Assume (for the sake of contradiction) that $\sigma\left(  P\right)
\subseteq\left[  n\right]  \setminus Q$. Thus,%
\begin{align*}
\left\vert \sigma\left(  P\right)  \right\vert  &  \leq\left\vert \left[
n\right]  \setminus Q\right\vert =\underbrace{\left\vert \left[  n\right]
\right\vert }_{=n}-\left\vert Q\right\vert \ \ \ \ \ \ \ \ \ \ \left(
\text{since }Q\subseteq\left[  n\right]  \right) \\
&  =n-\left\vert Q\right\vert <\left\vert P\right\vert
\ \ \ \ \ \ \ \ \ \ \left(  \text{since }\left\vert P\right\vert +\left\vert
Q\right\vert >n\right)  .
\end{align*}
But the map $\sigma$ is injective (since it is a permutation); thus,
$\left\vert \sigma\left(  P\right)  \right\vert =\left\vert P\right\vert $.
Hence, $\left\vert P\right\vert =\left\vert \sigma\left(  P\right)
\right\vert <\left\vert P\right\vert $, which is absurd. This contradiction
shows that our assumption (that $\sigma\left(  P\right)  \subseteq\left[
n\right]  \setminus Q$) was false. Hence, $\sigma\left(  P\right)
\not \subseteq \left[  n\right]  \setminus Q$. In other words, there exists
some $p\in P$ such that $\sigma\left(  p\right)  \notin\left[  n\right]
\setminus Q$. Consider this $p$.

From $\sigma\left(  p\right)  \in\left[  n\right]  $ and $\sigma\left(
p\right)  \notin\left[  n\right]  \setminus Q$, we obtain $\sigma\left(
p\right)  \in\left[  n\right]  \setminus\left(  \left[  n\right]  \setminus
Q\right)  \subseteq Q$. Hence, (\ref{eq.exe.det.blocktria-twisted.ass})
(applied to $i=p$ and $j=\sigma\left(  p\right)  $) yields $a_{p,\sigma\left(
p\right)  }=0$.

Now, $a_{p,\sigma\left(  p\right)  }$ is a factor of the product $\prod
_{i=1}^{n}a_{i,\sigma\left(  i\right)  }$. Thus, at least one factor of the
product $\prod_{i=1}^{n}a_{i,\sigma\left(  i\right)  }$ is $0$ (namely,
$a_{p,\sigma\left(  p\right)  }=0$). Hence, the whole product must be $0$. In
other words, $\prod_{i=1}^{n}a_{i,\sigma\left(  i\right)  }=0$.

Now, forget that we fixed $\sigma$. We thus have shown that $\prod_{i=1}%
^{n}a_{i,\sigma\left(  i\right)  }=0$ for each $\sigma\in S_{n}$. Hence,
(\ref{eq.det.eq.2}) becomes%
\[
\det A=\sum_{\sigma\in S_{n}}\left(  -1\right)  ^{\sigma}\underbrace{\prod
_{i=1}^{n}a_{i,\sigma\left(  i\right)  }}_{=0}=\sum_{\sigma\in S_{n}}\left(
-1\right)  ^{\sigma}0=0.
\]
This solves Exercise \ref{exe.det.blocktria-twisted} \textbf{(a)}.
\end{proof}

Let us now prepare for a \textquotedblleft proper\textquotedblright\ solution
to Exercise \ref{exe.det.blocktria-twisted}, which will solve both parts
\textbf{(a)} and \textbf{(b)}. We shall use the notations introduced in
Definition \ref{def.submatrix} and in Definition
\ref{def.sect.laplace.notations}. We shall also use the following lemma:

\begin{lemma}
\label{lem.sol.det.blocktria-twisted.0}Let $n\in\mathbb{N}$. For any subset
$I$ of $\left\{  1,2,\ldots,n\right\}  $, we let $\widetilde{I}$ denote the
complement $\left\{  1,2,\ldots,n\right\}  \setminus I$ of $I$. (For instance,
if $n=4$ and $I=\left\{  1,4\right\}  $, then $\widetilde{I}=\left\{
2,3\right\}  $.)

Let $P$ and $Q$ be two subsets of $\left\{  1,2,\ldots,n\right\}  $. Let $R$
be a subset of $\left\{  1,2,\ldots,n\right\}  $ satisfying $\left\vert
R\right\vert =\left\vert P\right\vert $ and $R\not \subseteq \widetilde{Q}$.

Let $A=\left(  a_{i,j}\right)  _{1\leq i\leq n,\ 1\leq j\leq n}\in
\mathbb{K}^{n\times n}$ be an $n\times n$-matrix such that%
\begin{equation}
\text{every }i\in P\text{ and }j\in Q\text{ satisfy }a_{i,j}=0.
\label{eq.lem.sol.det.blocktria-twisted.0.ass}%
\end{equation}
Then, $\det\left(  \operatorname*{sub}\nolimits_{w\left(  P\right)
}^{w\left(  R\right)  }A\right)  =0$.
\end{lemma}

\begin{proof}
[Proof of Lemma \ref{lem.sol.det.blocktria-twisted.0}.]Define $k\in\mathbb{N}$
by $k=\left\vert R\right\vert =\left\vert P\right\vert $. (This is
well-defined, since $\left\vert R\right\vert =\left\vert P\right\vert $.)

\begin{vershort}
Recall that $w\left(  P\right)  $ is the list of all elements of $P$ in
increasing order (with no repetitions) (by the definition of $w\left(
P\right)  $). Thus, $w\left(  P\right)  $ is a list of $k$ elements (since
$\left\vert P\right\vert =k$). Write this list $w\left(  P\right)  $ in the
form $w\left(  P\right)  =\left(  p_{1},p_{2},\ldots,p_{k}\right)  $. Thus,
$P=\left\{  p_{1},p_{2},\ldots,p_{k}\right\}  $.

Similarly, write the list $w\left(  R\right)  $ in the form $w\left(
R\right)  =\left(  r_{1},r_{2},\ldots,r_{k}\right)  $. Thus, $R=\left\{
r_{1},r_{2},\ldots,r_{k}\right\}  $.
\end{vershort}

\begin{verlong}
We know that $w\left(  P\right)  $ is the list of all elements of $P$ in
increasing order (with no repetitions) (by the definition of $w\left(
P\right)  $). Thus, $w\left(  P\right)  $ is a list of $\left\vert
P\right\vert $ elements. In other words, $w\left(  P\right)  $ is a list of
$k$ elements (since $\left\vert P\right\vert =k$).

Write $w\left(  P\right)  $ in the form $w\left(  P\right)  =\left(
p_{1},p_{2},\ldots,p_{k}\right)  $. (This is possible, since $w\left(
P\right)  $ is a list of $k$ elements.)

We know that $w\left(  P\right)  $ is a list of all elements of $P$. In other
words, $\left(  p_{1},p_{2},\ldots,p_{k}\right)  $ is a list of all elements
of $P$ (since $w\left(  P\right)  =\left(  p_{1},p_{2},\ldots,p_{k}\right)
$). Thus, $P=\left\{  p_{1},p_{2},\ldots,p_{k}\right\}  $.

We know that $w\left(  R\right)  $ is the list of all elements of $R$ in
increasing order (with no repetitions) (by the definition of $w\left(
R\right)  $). Thus, $w\left(  R\right)  $ is a list of $\left\vert
R\right\vert $ elements. In other words, $w\left(  R\right)  $ is a list of
$k$ elements (since $\left\vert R\right\vert =k$).

Write $w\left(  R\right)  $ in the form $w\left(  R\right)  =\left(
r_{1},r_{2},\ldots,r_{k}\right)  $. (This is possible, since $w\left(
R\right)  $ is a list of $k$ elements.)

We know that $w\left(  R\right)  $ is a list of all elements of $R$. In other
words, $\left(  r_{1},r_{2},\ldots,r_{k}\right)  $ is a list of all elements
of $R$ (since $w\left(  R\right)  =\left(  r_{1},r_{2},\ldots,r_{k}\right)
$). Thus, $R=\left\{  r_{1},r_{2},\ldots,r_{k}\right\}  $.
\end{verlong}

We have $R\not \subseteq \widetilde{Q}$. In other words, there exists some
$z\in R$ such that $z\notin\widetilde{Q}$. Consider this $z$. We have
$\widetilde{Q}=\left\{  1,2,\ldots,n\right\}  \setminus Q$ (by the definition
of $\widetilde{Q}$). Combining $z\in R\subseteq\left\{  1,2,\ldots,n\right\}
$ with $z\notin\widetilde{Q}$, we obtain%
\[
z\in\left\{  1,2,\ldots,n\right\}  \setminus\underbrace{\widetilde{Q}%
}_{=\left\{  1,2,\ldots,n\right\}  \setminus Q}=\left\{  1,2,\ldots,n\right\}
\setminus\left(  \left\{  1,2,\ldots,n\right\}  \setminus Q\right)  \subseteq
Q.
\]
But $z\in R=\left\{  r_{1},r_{2},\ldots,r_{k}\right\}  $. Hence, $z=r_{v}$ for
some $v\in\left\{  1,2,\ldots,k\right\}  $. Consider this $v$. We have%
\begin{equation}
a_{p_{x},r_{v}}=0\ \ \ \ \ \ \ \ \ \ \text{for each }x\in\left\{
1,2,\ldots,k\right\}  . \label{pf.lem.sol.det.blocktria-twisted.0.2}%
\end{equation}

[\textit{Proof of (\ref{pf.lem.sol.det.blocktria-twisted.0.2}):} Let
$x\in\left\{  1,2,\ldots,k\right\}  $. Then, $p_{x}\in\left\{  p_{1}%
,p_{2},\ldots,p_{k}\right\}  =P$. Also, $r_{v}=z\in Q$. Hence,
(\ref{eq.lem.sol.det.blocktria-twisted.0.ass}) (applied to $i=p_{x}$ and
$j=r_{v}$) yields $a_{p_{x},r_{v}}=0$. This proves
(\ref{pf.lem.sol.det.blocktria-twisted.0.2}).]

\begin{vershort}
From $w\left(  P\right)  =\left(  p_{1},p_{2},\ldots,p_{k}\right)  $ and
$w\left(  R\right)  =\left(  r_{1},r_{2},\ldots,r_{k}\right)  $, we obtain%
\[
\operatorname*{sub}\nolimits_{w\left(  P\right)  }^{w\left(  R\right)
}A=\operatorname*{sub}\nolimits_{\left(  p_{1},p_{2},\ldots,p_{k}\right)
}^{\left(  r_{1},r_{2},\ldots,r_{k}\right)  }A=\operatorname*{sub}%
\nolimits_{p_{1},p_{2},\ldots,p_{k}}^{r_{1},r_{2},\ldots,r_{k}}A=\left(
a_{p_{x},r_{y}}\right)  _{1\leq x\leq k,\ 1\leq y\leq k}.
\]
Hence, $\operatorname*{sub}\nolimits_{w\left(  P\right)  }^{w\left(  R\right)
}A$ is a $k\times k$-matrix, and its $v$-th column is%
\[
\left(
\begin{array}
[c]{c}%
a_{p_{1},r_{v}}\\
a_{p_{2},r_{v}}\\
\vdots\\
a_{p_{k},r_{v}}%
\end{array}
\right)  =\left(  \underbrace{a_{p_{x},r_{v}}}_{\substack{=0\\\text{(by
(\ref{pf.lem.sol.det.blocktria-twisted.0.2}))}}}\right)  _{1\leq x\leq
k,\ 1\leq y\leq1}=\left(  0\right)  _{1\leq x\leq k,\ 1\leq y\leq1}.
\]
In other words, the $v$-th column of the matrix $\operatorname*{sub}%
\nolimits_{w\left(  P\right)  }^{w\left(  R\right)  }A$ consists of zeroes.
Thus, a column of the matrix $\operatorname*{sub}\nolimits_{w\left(  P\right)
}^{w\left(  R\right)  }A$ consists of zeroes.
\end{vershort}

\begin{verlong}
From $w\left(  P\right)  =\left(  p_{1},p_{2},\ldots,p_{k}\right)  $ and
$w\left(  R\right)  =\left(  r_{1},r_{2},\ldots,r_{k}\right)  $, we obtain%
\begin{align}
\operatorname*{sub}\nolimits_{w\left(  P\right)  }^{w\left(  R\right)  }A  &
=\operatorname*{sub}\nolimits_{\left(  p_{1},p_{2},\ldots,p_{k}\right)
}^{\left(  r_{1},r_{2},\ldots,r_{k}\right)  }A=\operatorname*{sub}%
\nolimits_{p_{1},p_{2},\ldots,p_{k}}^{r_{1},r_{2},\ldots,r_{k}}%
A\ \ \ \ \ \ \ \ \ \ \left(  \text{by the definition of }\operatorname*{sub}%
\nolimits_{\left(  p_{1},p_{2},\ldots,p_{k}\right)  }^{\left(  r_{1}%
,r_{2},\ldots,r_{k}\right)  }A\right) \nonumber\\
&  =\left(  a_{p_{x},r_{y}}\right)  _{1\leq x\leq k,\ 1\leq y\leq k}
\label{pf.lem.sol.det.blocktria-twisted.0.3}%
\end{align}
(by the definition of $\operatorname*{sub}\nolimits_{p_{1},p_{2},\ldots,p_{k}%
}^{r_{1},r_{2},\ldots,r_{k}}A$, since $A=\left(  a_{i,j}\right)  _{1\leq i\leq
n,\ 1\leq j\leq n}$). Thus, $\operatorname*{sub}\nolimits_{w\left(  P\right)
}^{w\left(  R\right)  }A$ is a $k\times k$-matrix. Moreover, from
(\ref{pf.lem.sol.det.blocktria-twisted.0.3}), we obtain
\begin{align*}
\left(  \text{the }v\text{-th column of the matrix }\operatorname*{sub}%
\nolimits_{w\left(  P\right)  }^{w\left(  R\right)  }A\right)   &  =\left(
\begin{array}
[c]{c}%
a_{p_{1},r_{v}}\\
a_{p_{2},r_{v}}\\
\vdots\\
a_{p_{k},r_{v}}%
\end{array}
\right)  =\left(  \underbrace{a_{p_{x},r_{v}}}_{\substack{=0\\\text{(by
(\ref{pf.lem.sol.det.blocktria-twisted.0.2}))}}}\right)  _{1\leq x\leq
k,\ 1\leq y\leq1}\\
&  =\left(  0\right)  _{1\leq x\leq k,\ 1\leq y\leq1}.
\end{align*}
In other words, the $v$-th column of the matrix $\operatorname*{sub}%
\nolimits_{w\left(  P\right)  }^{w\left(  R\right)  }A$ consists of zeroes.
Thus, a column of the matrix $\operatorname*{sub}\nolimits_{w\left(  P\right)
}^{w\left(  R\right)  }A$ consists of zeroes.
\end{verlong}

Hence, Exercise \ref{exe.ps4.6} \textbf{(d)} (applied to $k$ and
$\operatorname*{sub}\nolimits_{w\left(  P\right)  }^{w\left(  R\right)  }A$
instead of $n$ and $A$) yields $\det\left(  \operatorname*{sub}%
\nolimits_{w\left(  P\right)  }^{w\left(  R\right)  }A\right)  =0$. This
proves Lemma \ref{lem.sol.det.blocktria-twisted.0}.
\end{proof}

\begin{proof}
[Solution to Exercise \ref{exe.det.blocktria-twisted}.]For any subset $I$ of
$\left\{  1,2,\ldots,n\right\}  $, we let $\widetilde{I}$ denote the
complement $\left\{  1,2,\ldots,n\right\}  \setminus I$ of $I$. (For instance,
if $n=4$ and $I=\left\{  1,4\right\}  $, then $\widetilde{I}=\left\{
2,3\right\}  $.)

\begin{vershort}
Thus, $\widetilde{Q}$ is the complement of $Q$. Hence, $Q$ is, in turn, the
complement of $\widetilde{Q}$. In other words, $Q=\widetilde{\widetilde{Q}}$.
Also, $\widetilde{Q}$ is the complement of $Q$ in the $n$-element set
$\left\{  1,2,\ldots,n\right\}  $; therefore, $\left\vert \widetilde{Q}%
\right\vert =n-\left\vert Q\right\vert $.
\end{vershort}

\begin{verlong}
The definition of $\widetilde{Q}$ yields $\widetilde{Q}=\left\{
1,2,\ldots,n\right\}  \setminus Q$, and thus%
\begin{align}
\left\vert \widetilde{Q}\right\vert  &  =\left\vert \left\{  1,2,\ldots
,n\right\}  \setminus Q\right\vert =\underbrace{\left\vert \left\{
1,2,\ldots,n\right\}  \right\vert }_{=n}-\left\vert Q\right\vert
\ \ \ \ \ \ \ \ \ \ \left(  \text{since }Q\text{ is a subset of }\left\{
1,2,\ldots,n\right\}  \right) \nonumber\\
&  =n-\left\vert Q\right\vert . \label{sol.det.blocktria-twisted.1}%
\end{align}
Also, $\widetilde{Q}=\left\{  1,2,\ldots,n\right\}  \setminus Q\subseteq
\left\{  1,2,\ldots,n\right\}  $. In other words, $\widetilde{Q}$ is a subset
of $\left\{  1,2,\ldots,n\right\}  $. Its complement $\widetilde{\widetilde{Q}%
}$ is%
\begin{align*}
\widetilde{\widetilde{Q}}  &  =\left\{  1,2,\ldots,n\right\}  \setminus
\underbrace{\widetilde{Q}}_{=\left\{  1,2,\ldots,n\right\}  \setminus
Q}\ \ \ \ \ \ \ \ \ \ \left(  \text{by the definition of }%
\widetilde{\widetilde{Q}}\right) \\
&  =\left\{  1,2,\ldots,n\right\}  \setminus\left(  \left\{  1,2,\ldots
,n\right\}  \setminus Q\right)  =Q
\end{align*}
(since $Q$ is a subset of $\left\{  1,2,\ldots,n\right\}  $).
\end{verlong}

For every subset $U$ of $\left\{  1,2,\ldots,n\right\}  $, we have
\begin{equation}
\det A=\sum_{\substack{R\subseteq\left\{  1,2,\ldots,n\right\}  ;\\\left\vert
R\right\vert =\left\vert U\right\vert }}\left(  -1\right)  ^{\sum U+\sum
R}\det\left(  \operatorname*{sub}\nolimits_{w\left(  U\right)  }^{w\left(
R\right)  }A\right)  \det\left(  \operatorname*{sub}\nolimits_{w\left(
\widetilde{U}\right)  }^{w\left(  \widetilde{R}\right)  }A\right)  .
\label{sol.det.blocktria-twisted.lap}%
\end{equation}
(Indeed, this is merely the claim of Theorem \ref{thm.det.laplace-multi}
\textbf{(a)}, with the variables $P$ and $Q$ renamed as $U$ and $R$.)

\textbf{(a)} Assume that $\left\vert P\right\vert +\left\vert Q\right\vert
>n$. Thus, $\left\vert P\right\vert >n-\left\vert Q\right\vert =\left\vert
\widetilde{Q}\right\vert $.

Hence, each subset $R$ of $\left\{  1,2,\ldots,n\right\}  $ satisfying
$\left\vert R\right\vert =\left\vert P\right\vert $ must satisfy%
\begin{equation}
\det\left(  \operatorname*{sub}\nolimits_{w\left(  P\right)  }^{w\left(
R\right)  }A\right)  =0. \label{sol.det.blocktria-twisted.a.1}%
\end{equation}

[\textit{Proof of (\ref{sol.det.blocktria-twisted.a.1}):} Let $R$ be a subset
of $\left\{  1,2,\ldots,n\right\}  $ satisfying $\left\vert R\right\vert
=\left\vert P\right\vert $. Then, $\left\vert R\right\vert =\left\vert
P\right\vert >\left\vert \widetilde{Q}\right\vert $. If we had $R\subseteq
\widetilde{Q}$, then we would have $\left\vert R\right\vert \leq\left\vert
\widetilde{Q}\right\vert $, which would contradict $\left\vert R\right\vert
>\left\vert \widetilde{Q}\right\vert $. Thus, we cannot have $R\subseteq
\widetilde{Q}$. Hence, we have $R\not \subseteq \widetilde{Q}$. Hence, Lemma
\ref{lem.sol.det.blocktria-twisted.0} yields $\det\left(  \operatorname*{sub}%
\nolimits_{w\left(  P\right)  }^{w\left(  R\right)  }A\right)  =0$. This
proves (\ref{sol.det.blocktria-twisted.a.1}).]

Now, (\ref{sol.det.blocktria-twisted.lap}) (applied to $U=P$) yields%
\begin{align*}
\det A  &  =\sum_{\substack{R\subseteq\left\{  1,2,\ldots,n\right\}
;\\\left\vert R\right\vert =\left\vert P\right\vert }}\left(  -1\right)
^{\sum P+\sum R}\underbrace{\det\left(  \operatorname*{sub}\nolimits_{w\left(
P\right)  }^{w\left(  R\right)  }A\right)  }_{\substack{=0\\\text{(by
(\ref{sol.det.blocktria-twisted.a.1}))}}}\det\left(  \operatorname*{sub}%
\nolimits_{w\left(  \widetilde{P}\right)  }^{w\left(  \widetilde{R}\right)
}A\right) \\
&  =\sum_{\substack{R\subseteq\left\{  1,2,\ldots,n\right\}  ;\\\left\vert
R\right\vert =\left\vert P\right\vert }}\left(  -1\right)  ^{\sum P+\sum
R}0\det\left(  \operatorname*{sub}\nolimits_{w\left(  \widetilde{P}\right)
}^{w\left(  \widetilde{R}\right)  }A\right)  =0.
\end{align*}
This solves Exercise \ref{exe.det.blocktria-twisted} \textbf{(a)}.

\textbf{(b)} Assume that $\left\vert P\right\vert +\left\vert Q\right\vert
=n$. Thus, $\left\vert P\right\vert =n-\left\vert Q\right\vert =\left\vert
\widetilde{Q}\right\vert $. In other words, $\left\vert \widetilde{Q}%
\right\vert =\left\vert P\right\vert $. Hence, $\widetilde{Q}$ is a subset $R$
of $\left\{  1,2,\ldots,n\right\}  $ satisfying $\left\vert R\right\vert
=\left\vert P\right\vert $ (since $\widetilde{Q}$ is a subset of $\left\{
1,2,\ldots,n\right\}  $).

Hence, each subset $R$ of $\left\{  1,2,\ldots,n\right\}  $ satisfying
$\left\vert R\right\vert =\left\vert P\right\vert $ and $R\neq\widetilde{Q}$
must satisfy%
\begin{equation}
\det\left(  \operatorname*{sub}\nolimits_{w\left(  P\right)  }^{w\left(
R\right)  }A\right)  =0. \label{sol.det.blocktria-twisted.b.1}%
\end{equation}

[\textit{Proof of (\ref{sol.det.blocktria-twisted.b.1}):} Let $R$ be a subset
of $\left\{  1,2,\ldots,n\right\}  $ satisfying $\left\vert R\right\vert
=\left\vert P\right\vert $ and $R\neq\widetilde{Q}$. Thus, $\left\vert
R\right\vert =\left\vert P\right\vert =\left\vert \widetilde{Q}\right\vert $.

\begin{vershort}
Thus, the set $R$ has the same size as $\widetilde{Q}$. If $R$ was a subset of
$\widetilde{Q}$, then this would lead to $R=\widetilde{Q}$ (because the only
subset of $\widetilde{Q}$ having the same size as $\widetilde{Q}$ is
$\widetilde{Q}$ itself), which would contradict $R\neq\widetilde{Q}$. Hence,
$R$ is not a subset of $\widetilde{Q}$.
\end{vershort}

\begin{verlong}
Assume (for the sake of contradiction) that $R$ is a subset of $\widetilde{Q}$.

But the following fact is well-known: If $X$ is a finite set, and if $Y$ is a
subset of $X$ such that $\left\vert Y\right\vert \geq\left\vert X\right\vert
$, then $Y=X$. Applying this to $X=\widetilde{Q}$ and $Y=R$, we conclude that
$R=\widetilde{Q}$ (since $R$ is a subset of $\widetilde{Q}$ and since
$\left\vert R\right\vert =\left\vert \widetilde{Q}\right\vert \geq\left\vert
\widetilde{Q}\right\vert $). This contradicts $R\neq\widetilde{Q}$. This
contradiction shows that our assumption (that $R$ is a subset of
$\widetilde{Q}$) was false. Hence, $R$ is not a subset of $\widetilde{Q}$.
\end{verlong}

In other words, $R\not \subseteq \widetilde{Q}$. Hence, Lemma
\ref{lem.sol.det.blocktria-twisted.0} yields $\det\left(  \operatorname*{sub}%
\nolimits_{w\left(  P\right)  }^{w\left(  R\right)  }A\right)  =0$. This
proves (\ref{sol.det.blocktria-twisted.b.1}).]

Now, (\ref{sol.det.blocktria-twisted.lap}) (applied to $U=P$) yields%
\begin{align*}
\det A  &  =\sum_{\substack{R\subseteq\left\{  1,2,\ldots,n\right\}
;\\\left\vert R\right\vert =\left\vert P\right\vert }}\left(  -1\right)
^{\sum P+\sum R}\det\left(  \operatorname*{sub}\nolimits_{w\left(  P\right)
}^{w\left(  R\right)  }A\right)  \det\left(  \operatorname*{sub}%
\nolimits_{w\left(  \widetilde{P}\right)  }^{w\left(  \widetilde{R}\right)
}A\right) \\
&  =\left(  -1\right)  ^{\sum P+\sum\widetilde{Q}}\det\left(
\operatorname*{sub}\nolimits_{w\left(  P\right)  }^{w\left(  \widetilde{Q}%
\right)  }A\right)  \underbrace{\det\left(  \operatorname*{sub}%
\nolimits_{w\left(  \widetilde{P}\right)  }^{w\left(  \widetilde{\widetilde{Q}%
}\right)  }A\right)  }_{\substack{=\det\left(  \operatorname*{sub}%
\nolimits_{w\left(  \widetilde{P}\right)  }^{w\left(  Q\right)  }A\right)
\\\text{(since }\widetilde{\widetilde{Q}}=Q\text{)}}}\\
&  \ \ \ \ \ \ \ \ \ \ +\sum_{\substack{R\subseteq\left\{  1,2,\ldots
,n\right\}  ;\\\left\vert R\right\vert =\left\vert P\right\vert ;\\R\neq
\widetilde{Q}}}\left(  -1\right)  ^{\sum P+\sum R}\underbrace{\det\left(
\operatorname*{sub}\nolimits_{w\left(  P\right)  }^{w\left(  R\right)
}A\right)  }_{\substack{=0\\\text{(by (\ref{sol.det.blocktria-twisted.b.1}))}%
}}\det\left(  \operatorname*{sub}\nolimits_{w\left(  \widetilde{P}\right)
}^{w\left(  \widetilde{R}\right)  }A\right) \\
&  \ \ \ \ \ \ \ \ \ \ \left(
\begin{array}
[c]{c}%
\text{here, we have split off the addend for }R=\widetilde{Q}\text{ from the
sum}\\
\text{(since }\widetilde{Q}\text{ is a subset }R\text{ of }\left\{
1,2,\ldots,n\right\}  \text{ satisfying }\left\vert R\right\vert =\left\vert
P\right\vert \text{)}%
\end{array}
\right) \\
&  =\left(  -1\right)  ^{\sum P+\sum\widetilde{Q}}\det\left(
\operatorname*{sub}\nolimits_{w\left(  P\right)  }^{w\left(  \widetilde{Q}%
\right)  }A\right)  \det\left(  \operatorname*{sub}\nolimits_{w\left(
\widetilde{P}\right)  }^{w\left(  Q\right)  }A\right) \\
&  \ \ \ \ \ \ \ \ \ \ +\underbrace{\sum_{\substack{R\subseteq\left\{
1,2,\ldots,n\right\}  ;\\\left\vert R\right\vert =\left\vert P\right\vert
;\\R\neq\widetilde{Q}}}\left(  -1\right)  ^{\sum P+\sum R}0\det\left(
\operatorname*{sub}\nolimits_{w\left(  \widetilde{P}\right)  }^{w\left(
\widetilde{R}\right)  }A\right)  }_{=0}\\
&  =\left(  -1\right)  ^{\sum P+\sum\widetilde{Q}}\det\left(
\operatorname*{sub}\nolimits_{w\left(  P\right)  }^{w\left(  \widetilde{Q}%
\right)  }A\right)  \det\left(  \operatorname*{sub}\nolimits_{w\left(
\widetilde{P}\right)  }^{w\left(  Q\right)  }A\right)  .
\end{align*}
This solves Exercise \ref{exe.det.blocktria-twisted} \textbf{(b)}.
\end{proof}

\subsection{Solution to Exercise \ref{exe.det(A+B)}}

\begin{vershort}
\begin{proof}
[Proof of Theorem \ref{thm.det(A+B)}.]Write the $n\times n$-matrix $A$ in the
form $A=\left(  a_{i,j}\right)  _{1\leq i\leq n,\ 1\leq j\leq n}$. Write the
$n\times n$-matrix $B$ in the form $B=\left(  b_{i,j}\right)  _{1\leq i\leq
n,\ 1\leq j\leq n}$.

Let $\left[  n\right]  $ denote the set $\left\{  1,2,\ldots,n\right\}  $.
Then, every subset $I$ of $\left[  n\right]  $ satisfies
\begin{equation}
\left[  n\right]  \setminus I=\widetilde{I} \label{sol.det(A+B).short.compl}%
\end{equation}
(because the definition of $\widetilde{I}$ yields $\widetilde{I}%
=\underbrace{\left\{  1,2,\ldots,n\right\}  }_{=\left[  n\right]  }\setminus
I=\left[  n\right]  \setminus I$).

For every $\sigma\in S_{n}$ and every subset $P$ of $\left\{  1,2,\ldots
,n\right\}  $, define an element $c_{\sigma,P}$ of $\mathbb{K}$ by
\begin{equation}
c_{\sigma,P}=\left(  -1\right)  ^{\sigma}\left(  \prod_{i\in P}a_{i,\sigma
\left(  i\right)  }\right)  \left(  \prod_{i\in\widetilde{P}}b_{i,\sigma
\left(  i\right)  }\right)  . \label{sol.det(A+B).short.c=}%
\end{equation}

If $P$ and $Q$ are two subsets of $\left\{  1,2,\ldots,n\right\}  $ satisfying
$\left\vert Q\right\vert \neq\left\vert P\right\vert $, then%
\begin{equation}
\sum_{\substack{\sigma\in S_{n};\\\sigma\left(  P\right)  =Q}}c_{\sigma,P}=0
\label{sol.det(A+B).short.emptysum}%
\end{equation}
\footnote{\textit{Proof of (\ref{sol.det(A+B).short.emptysum}):} This can be
shown in precisely the same way as (\ref{pf.thm.det.laplace-multi.emptysum})
was shown in our proof of Theorem \ref{thm.det.laplace-multi}.}.

If $P$ and $Q$ are two subsets of $\left\{  1,2,\ldots,n\right\}  $ satisfying
$\left\vert Q\right\vert =\left\vert P\right\vert $, then%
\begin{equation}
\sum_{\substack{\sigma\in S_{n};\\\sigma\left(  P\right)  =Q}}c_{\sigma
,P}=\left(  -1\right)  ^{\sum P+\sum Q}\det\left(  \operatorname*{sub}%
\nolimits_{w\left(  P\right)  }^{w\left(  Q\right)  }A\right)  \det\left(
\operatorname*{sub}\nolimits_{w\left(  \widetilde{P}\right)  }^{w\left(
\widetilde{Q}\right)  }B\right)  \label{sol.det(A+B).short.nonemptysum}%
\end{equation}
\footnote{\textit{Proof of (\ref{sol.det(A+B).short.nonemptysum}):} Let $P$
and $Q$ be two subsets of $\left\{  1,2,\ldots,n\right\}  $ satisfying
$\left\vert Q\right\vert =\left\vert P\right\vert $. From $\left\vert
Q\right\vert =\left\vert P\right\vert $, we obtain $\left\vert P\right\vert
=\left\vert Q\right\vert $. Lemma \ref{lem.det.laplace-multi.Apq} thus yields%
\[
\sum_{\substack{\sigma\in S_{n};\\\sigma\left(  P\right)  =Q}}\left(
-1\right)  ^{\sigma}\left(  \prod_{i\in P}a_{i,\sigma\left(  i\right)
}\right)  \left(  \prod_{i\in\widetilde{P}}b_{i,\sigma\left(  i\right)
}\right)  =\left(  -1\right)  ^{\sum P+\sum Q}\det\left(  \operatorname*{sub}%
\nolimits_{w\left(  P\right)  }^{w\left(  Q\right)  }A\right)  \det\left(
\operatorname*{sub}\nolimits_{w\left(  \widetilde{P}\right)  }^{w\left(
\widetilde{Q}\right)  }B\right)  .
\]
But this is precisely the equality (\ref{sol.det(A+B).short.nonemptysum})
(because of (\ref{sol.det(A+B).short.c=})). Thus,
(\ref{sol.det(A+B).short.nonemptysum}) is proven.}.

Adding the equalities $A=\left(  a_{i,j}\right)  _{1\leq i\leq n,\ 1\leq j\leq
n}$ and $B=\left(  b_{i,j}\right)  _{1\leq i\leq n,\ 1\leq j\leq n}$, we
obtain $A+B=\left(  a_{i,j}+b_{i,j}\right)  _{1\leq i\leq n,\ 1\leq j\leq n}$.
Thus, (\ref{eq.det.eq.2}) (applied to $A+B$ and $a_{i,j}+b_{i,j}$ instead of
$A$ and $a_{i,j}$) yields%
\begin{align}
\det\left(  A+B\right)   &  =\sum_{\sigma\in S_{n}}\left(  -1\right)
^{\sigma}\underbrace{\prod_{i=1}^{n}\left(  a_{i,\sigma\left(  i\right)
}+b_{i,\sigma\left(  i\right)  }\right)  }_{\substack{=\sum_{I\subseteq\left[
n\right]  }\left(  \prod_{i\in I}a_{i,\sigma\left(  i\right)  }\right)
\left(  \prod_{i\in\left[  n\right]  \setminus I}b_{i,\sigma\left(  i\right)
}\right)  \\\text{(by Exercise \ref{exe.prod(ai+bi)} \textbf{(a)}%
,}\\\text{applied to }a_{i,\sigma\left(  i\right)  }\text{ and }%
b_{i,\sigma\left(  i\right)  }\text{ instead of }a_{i}\text{ and }%
b_{i}\text{)}}}\nonumber\\
&  =\sum_{\sigma\in S_{n}}\left(  -1\right)  ^{\sigma}\underbrace{\sum
_{I\subseteq\left[  n\right]  }}_{=\sum_{I\subseteq\left\{  1,2,\ldots
,n\right\}  }}\left(  \prod_{i\in I}a_{i,\sigma\left(  i\right)  }\right)
\left(  \underbrace{\prod_{i\in\left[  n\right]  \setminus I}}%
_{\substack{=\prod_{i\in\widetilde{I}}\\\text{(by
(\ref{sol.det(A+B).short.compl}))}}}b_{i,\sigma\left(  i\right)  }\right)
\nonumber\\
&  =\sum_{\sigma\in S_{n}}\left(  -1\right)  ^{\sigma}\underbrace{\sum
_{I\subseteq\left\{  1,2,\ldots,n\right\}  }\left(  \prod_{i\in I}%
a_{i,\sigma\left(  i\right)  }\right)  \left(  \prod_{i\in\widetilde{I}%
}b_{i,\sigma\left(  i\right)  }\right)  }_{\substack{=\sum_{P\subseteq\left\{
1,2,\ldots,n\right\}  }\left(  \prod_{i\in P}a_{i,\sigma\left(  i\right)
}\right)  \left(  \prod_{i\in\widetilde{P}}b_{i,\sigma\left(  i\right)
}\right)  \\\text{(here, we have renamed the summation index }I\text{ as
}P\text{)}}}\nonumber\\
&  =\sum_{\sigma\in S_{n}}\left(  -1\right)  ^{\sigma}\sum_{P\subseteq\left\{
1,2,\ldots,n\right\}  }\left(  \prod_{i\in P}a_{i,\sigma\left(  i\right)
}\right)  \left(  \prod_{i\in\widetilde{P}}b_{i,\sigma\left(  i\right)
}\right) \nonumber\\
&  =\sum_{P\subseteq\left\{  1,2,\ldots,n\right\}  }\sum_{\sigma\in S_{n}%
}\underbrace{\left(  -1\right)  ^{\sigma}\left(  \prod_{i\in P}a_{i,\sigma
\left(  i\right)  }\right)  \left(  \prod_{i\in\widetilde{P}}b_{i,\sigma
\left(  i\right)  }\right)  }_{\substack{=c_{\sigma,P}\\\text{(by
(\ref{sol.det(A+B).short.c=}))}}}\nonumber\\
&  =\sum_{P\subseteq\left\{  1,2,\ldots,n\right\}  }\sum_{\sigma\in S_{n}%
}c_{\sigma,P}. \label{sol.det(A+B).short.1}%
\end{align}

But every subset $P$ of $\left\{  1,2,\ldots,n\right\}  $ satisfies%
\begin{align*}
&  \underbrace{\sum_{\sigma\in S_{n}}}_{\substack{=\sum_{Q\subseteq\left\{
1,2,\ldots,n\right\}  }\sum_{\substack{\sigma\in S_{n};\\\sigma\left(
P\right)  =Q}}\\\text{(because for every }\sigma\in S_{n}\text{, the
set}\\\sigma\left(  P\right)  \text{ is a subset of }\left\{  1,2,\ldots
,n\right\}  \text{)}}}c_{\sigma,P}\\
&  =\sum_{Q\subseteq\left\{  1,2,\ldots,n\right\}  }\sum_{\substack{\sigma\in
S_{n};\\\sigma\left(  P\right)  =Q}}c_{\sigma,P}=\sum_{\substack{Q\subseteq
\left\{  1,2,\ldots,n\right\}  ;\\\left\vert Q\right\vert =\left\vert
P\right\vert }}\sum_{\substack{\sigma\in S_{n};\\\sigma\left(  P\right)
=Q}}c_{\sigma,P}+\sum_{\substack{Q\subseteq\left\{  1,2,\ldots,n\right\}
;\\\left\vert Q\right\vert \neq\left\vert P\right\vert }}\underbrace{\sum
_{\substack{\sigma\in S_{n};\\\sigma\left(  P\right)  =Q}}c_{\sigma,P}%
}_{\substack{=0\\\text{(by (\ref{sol.det(A+B).short.emptysum}))}}}\\
&  \ \ \ \ \ \ \ \ \ \ \left(
\begin{array}
[c]{c}%
\text{since every subset }Q\text{ of }\left\{  1,2,\ldots,n\right\}  \text{
satisfies}\\
\text{either }\left\vert Q\right\vert =\left\vert P\right\vert \text{ or
}\left\vert Q\right\vert \neq\left\vert P\right\vert \text{ (but not both)}%
\end{array}
\right) \\
&  =\sum_{\substack{Q\subseteq\left\{  1,2,\ldots,n\right\}  ;\\\left\vert
Q\right\vert =\left\vert P\right\vert }}\sum_{\substack{\sigma\in
S_{n};\\\sigma\left(  P\right)  =Q}}c_{\sigma,P}+\underbrace{\sum
_{\substack{Q\subseteq\left\{  1,2,\ldots,n\right\}  ;\\\left\vert
Q\right\vert \neq\left\vert P\right\vert }}0}_{=0}\\
&  =\sum_{\substack{Q\subseteq\left\{  1,2,\ldots,n\right\}  ;\\\left\vert
Q\right\vert =\left\vert P\right\vert }}\underbrace{\sum_{\substack{\sigma\in
S_{n};\\\sigma\left(  P\right)  =Q}}c_{\sigma,P}}_{\substack{=\left(
-1\right)  ^{\sum P+\sum Q}\det\left(  \operatorname*{sub}\nolimits_{w\left(
P\right)  }^{w\left(  Q\right)  }A\right)  \det\left(  \operatorname*{sub}%
\nolimits_{w\left(  \widetilde{P}\right)  }^{w\left(  \widetilde{Q}\right)
}B\right)  \\\text{(by (\ref{sol.det(A+B).short.nonemptysum}))}}}\\
&  =\sum_{\substack{Q\subseteq\left\{  1,2,\ldots,n\right\}  ;\\\left\vert
Q\right\vert =\left\vert P\right\vert }}\left(  -1\right)  ^{\sum P+\sum
Q}\det\left(  \operatorname*{sub}\nolimits_{w\left(  P\right)  }^{w\left(
Q\right)  }A\right)  \det\left(  \operatorname*{sub}\nolimits_{w\left(
\widetilde{P}\right)  }^{w\left(  \widetilde{Q}\right)  }B\right)  .
\end{align*}
Hence, (\ref{sol.det(A+B).short.1}) becomes%
\begin{align*}
\det\left(  A+B\right)   &  =\sum_{P\subseteq\left\{  1,2,\ldots,n\right\}
}\underbrace{\sum_{\sigma\in S_{n}}c_{\sigma,P}}_{=\sum_{\substack{Q\subseteq
\left\{  1,2,\ldots,n\right\}  ;\\\left\vert Q\right\vert =\left\vert
P\right\vert }}\left(  -1\right)  ^{\sum P+\sum Q}\det\left(
\operatorname*{sub}\nolimits_{w\left(  P\right)  }^{w\left(  Q\right)
}A\right)  \det\left(  \operatorname*{sub}\nolimits_{w\left(  \widetilde{P}%
\right)  }^{w\left(  \widetilde{Q}\right)  }B\right)  }\\
&  =\sum_{P\subseteq\left\{  1,2,\ldots,n\right\}  }\sum_{\substack{Q\subseteq
\left\{  1,2,\ldots,n\right\}  ;\\\left\vert Q\right\vert =\left\vert
P\right\vert }}\left(  -1\right)  ^{\sum P+\sum Q}\det\left(
\operatorname*{sub}\nolimits_{w\left(  P\right)  }^{w\left(  Q\right)
}A\right)  \det\left(  \operatorname*{sub}\nolimits_{w\left(  \widetilde{P}%
\right)  }^{w\left(  \widetilde{Q}\right)  }B\right)  .
\end{align*}
This proves Theorem \ref{thm.det(A+B)}.
\end{proof}
\end{vershort}

\begin{verlong}
\begin{proof}
[Proof of Theorem \ref{thm.det(A+B)}.]Write the $n\times n$-matrix $A$ in the
form $A=\left(  a_{i,j}\right)  _{1\leq i\leq n,\ 1\leq j\leq n}$.

Write the $n\times n$-matrix $B$ in the form $B=\left(  b_{i,j}\right)
_{1\leq i\leq n,\ 1\leq j\leq n}$.

Let $\left[  n\right]  $ denote the set $\left\{  1,2,\ldots,n\right\}  $.
Then,
\begin{equation}
\left[  n\right]  \setminus I=\widetilde{I}\ \ \ \ \ \ \ \ \ \ \text{for every
subset }I\text{ of }\left[  n\right]  \label{sol.det(A+B).compl}%
\end{equation}
\footnote{\textit{Proof of (\ref{sol.det(A+B).compl}):} Let $I$ be a subset of
$\left[  n\right]  $. Then, $I\subseteq\left[  n\right]  =\left\{
1,2,\ldots,n\right\}  $. Now, the definition of $\widetilde{I}$ yields
$\widetilde{I}=\underbrace{\left\{  1,2,\ldots,n\right\}  }_{=\left[
n\right]  }\setminus I=\left[  n\right]  \setminus I$. This proves
(\ref{sol.det(A+B).compl}).}.

For every $\sigma\in S_{n}$ and every subset $P$ of $\left\{  1,2,\ldots
,n\right\}  $, define an element $c_{\sigma,P}$ of $\mathbb{K}$ by
\begin{equation}
c_{\sigma,P}=\left(  -1\right)  ^{\sigma}\left(  \prod_{i\in P}a_{i,\sigma
\left(  i\right)  }\right)  \left(  \prod_{i\in\widetilde{P}}b_{i,\sigma
\left(  i\right)  }\right)  . \label{sol.det(A+B).c=}%
\end{equation}

If $P$ and $Q$ are two subsets of $\left\{  1,2,\ldots,n\right\}  $ satisfying
$\left\vert Q\right\vert \neq\left\vert P\right\vert $, then%
\begin{equation}
\sum_{\substack{\sigma\in S_{n};\\\sigma\left(  P\right)  =Q}}c_{\sigma,P}=0
\label{sol.det(A+B).emptysum}%
\end{equation}
\footnote{\textit{Proof of (\ref{sol.det(A+B).emptysum}):} Let $P$ and $Q$ be
two subsets of $\left\{  1,2,\ldots,n\right\}  $ satisfying $\left\vert
Q\right\vert \neq\left\vert P\right\vert $.
\par
Let $\sigma\in S_{n}$ be such that $\sigma\left(  P\right)  =Q$. We shall
derive a contradiction.
\par
Indeed, Lemma \ref{lem.sol.det.laplace-multi.4} \textbf{(a)} shows that the
set $\sigma\left(  P\right)  $ is a subset of $\left\{  1,2,\ldots,n\right\}
$ satisfying $\left\vert \sigma\left(  P\right)  \right\vert =\left\vert
P\right\vert $. Hence, $\left\vert P\right\vert =\left\vert \underbrace{\sigma
\left(  P\right)  }_{=Q}\right\vert =\left\vert Q\right\vert \neq\left\vert
P\right\vert $. This is absurd. Hence, we have found a contradiction.
\par
Now, forget that we fixed $\sigma$. We thus have found a contradiction for
every $\sigma\in S_{n}$ satisfying $\sigma\left(  P\right)  =Q$. Thus, there
exists no $\sigma\in S_{n}$ satisfying $\sigma\left(  P\right)  =Q$. Hence,
the sum $\sum_{\substack{\sigma\in S_{n};\\\sigma\left(  P\right)
=Q}}c_{\sigma,P}$ is an empty sum. Thus, $\sum_{\substack{\sigma\in
S_{n};\\\sigma\left(  P\right)  =Q}}c_{\sigma,P}=\left(  \text{empty
sum}\right)  =0$. This proves (\ref{sol.det(A+B).emptysum}).}.

If $P$ and $Q$ are two subsets of $\left\{  1,2,\ldots,n\right\}  $ satisfying
$\left\vert Q\right\vert =\left\vert P\right\vert $, then%
\begin{equation}
\sum_{\substack{\sigma\in S_{n};\\\sigma\left(  P\right)  =Q}}c_{\sigma
,P}=\left(  -1\right)  ^{\sum P+\sum Q}\det\left(  \operatorname*{sub}%
\nolimits_{w\left(  P\right)  }^{w\left(  Q\right)  }A\right)  \det\left(
\operatorname*{sub}\nolimits_{w\left(  \widetilde{P}\right)  }^{w\left(
\widetilde{Q}\right)  }B\right)  \label{sol.det(A+B).nonemptysum}%
\end{equation}
\footnote{\textit{Proof of (\ref{sol.det(A+B).nonemptysum}):} Let $P$ and $Q$
be two subsets of $\left\{  1,2,\ldots,n\right\}  $ satisfying $\left\vert
Q\right\vert =\left\vert P\right\vert $. From $\left\vert Q\right\vert
=\left\vert P\right\vert $, we obtain $\left\vert P\right\vert =\left\vert
Q\right\vert $. Lemma \ref{lem.det.laplace-multi.Apq} thus yields%
\begin{align*}
&  \sum_{\substack{\sigma\in S_{n};\\\sigma\left(  P\right)  =Q}}\left(
-1\right)  ^{\sigma}\left(  \prod_{i\in P}a_{i,\sigma\left(  i\right)
}\right)  \left(  \prod_{i\in\widetilde{P}}b_{i,\sigma\left(  i\right)
}\right) \\
&  =\left(  -1\right)  ^{\sum P+\sum Q}\det\left(  \operatorname*{sub}%
\nolimits_{w\left(  P\right)  }^{w\left(  Q\right)  }A\right)  \det\left(
\operatorname*{sub}\nolimits_{w\left(  \widetilde{P}\right)  }^{w\left(
\widetilde{Q}\right)  }B\right)  .
\end{align*}
Now,%
\begin{align*}
&  \sum_{\substack{\sigma\in S_{n};\\\sigma\left(  P\right)  =Q}%
}\underbrace{c_{\sigma,P}}_{=\left(  -1\right)  ^{\sigma}\left(  \prod_{i\in
P}a_{i,\sigma\left(  i\right)  }\right)  \left(  \prod_{i\in\widetilde{P}%
}b_{i,\sigma\left(  i\right)  }\right)  }\\
&  =\sum_{\substack{\sigma\in S_{n};\\\sigma\left(  P\right)  =Q}}\left(
-1\right)  ^{\sigma}\left(  \prod_{i\in P}a_{i,\sigma\left(  i\right)
}\right)  \left(  \prod_{i\in\widetilde{P}}b_{i,\sigma\left(  i\right)
}\right) \\
&  =\left(  -1\right)  ^{\sum P+\sum Q}\det\left(  \operatorname*{sub}%
\nolimits_{w\left(  P\right)  }^{w\left(  Q\right)  }A\right)  \det\left(
\operatorname*{sub}\nolimits_{w\left(  \widetilde{P}\right)  }^{w\left(
\widetilde{Q}\right)  }B\right)  .
\end{align*}
This proves (\ref{sol.det(A+B).nonemptysum}).}.

Adding the equalities $A=\left(  a_{i,j}\right)  _{1\leq i\leq n,\ 1\leq j\leq
n}$ and $B=\left(  b_{i,j}\right)  _{1\leq i\leq n,\ 1\leq j\leq n}$, we
obtain%
\[
A+B=\left(  a_{i,j}\right)  _{1\leq i\leq n,\ 1\leq j\leq n}+\left(
b_{i,j}\right)  _{1\leq i\leq n,\ 1\leq j\leq n}=\left(  a_{i,j}%
+b_{i,j}\right)  _{1\leq i\leq n,\ 1\leq j\leq n}%
\]
(by the definition of the sum of two matrices). Thus, (\ref{eq.det.eq.2})
(applied to $A+B$ and $a_{i,j}+b_{i,j}$ instead of $A$ and $a_{i,j}$) yields%
\begin{align}
\det\left(  A+B\right)   &  =\sum_{\sigma\in S_{n}}\left(  -1\right)
^{\sigma}\underbrace{\prod_{i=1}^{n}\left(  a_{i,\sigma\left(  i\right)
}+b_{i,\sigma\left(  i\right)  }\right)  }_{\substack{=\sum_{I\subseteq\left[
n\right]  }\left(  \prod_{i\in I}a_{i,\sigma\left(  i\right)  }\right)
\left(  \prod_{i\in\left[  n\right]  \setminus I}b_{i,\sigma\left(  i\right)
}\right)  \\\text{(by Exercise \ref{exe.prod(ai+bi)} \textbf{(a)}%
}\\\text{(applied to }a_{i,\sigma\left(  i\right)  }\text{ and }%
b_{i,\sigma\left(  i\right)  }\text{ instead of }a_{i}\text{ and }%
b_{i}\text{))}}}\nonumber\\
&  =\sum_{\sigma\in S_{n}}\left(  -1\right)  ^{\sigma}\underbrace{\sum
_{I\subseteq\left[  n\right]  }}_{\substack{=\sum_{I\subseteq\left\{
1,2,\ldots,n\right\}  }\\\text{(since }\left[  n\right]  =\left\{
1,2,\ldots,n\right\}  \text{)}}}\left(  \prod_{i\in I}a_{i,\sigma\left(
i\right)  }\right)  \left(  \underbrace{\prod_{i\in\left[  n\right]  \setminus
I}}_{\substack{=\prod_{i\in\widetilde{I}}\\\text{(since }\left[  n\right]
\setminus I=\widetilde{I}\\\text{(by (\ref{sol.det(A+B).compl})))}%
}}b_{i,\sigma\left(  i\right)  }\right) \nonumber\\
&  =\sum_{\sigma\in S_{n}}\left(  -1\right)  ^{\sigma}\underbrace{\sum
_{I\subseteq\left\{  1,2,\ldots,n\right\}  }\left(  \prod_{i\in I}%
a_{i,\sigma\left(  i\right)  }\right)  \left(  \prod_{i\in\widetilde{I}%
}b_{i,\sigma\left(  i\right)  }\right)  }_{\substack{=\sum_{P\subseteq\left\{
1,2,\ldots,n\right\}  }\left(  \prod_{i\in P}a_{i,\sigma\left(  i\right)
}\right)  \left(  \prod_{i\in\widetilde{P}}b_{i,\sigma\left(  i\right)
}\right)  \\\text{(here, we have renamed the summation index }I\text{ as
}P\text{)}}}\nonumber\\
&  =\sum_{\sigma\in S_{n}}\underbrace{\left(  -1\right)  ^{\sigma}%
\sum_{P\subseteq\left\{  1,2,\ldots,n\right\}  }\left(  \prod_{i\in
P}a_{i,\sigma\left(  i\right)  }\right)  \left(  \prod_{i\in\widetilde{P}%
}b_{i,\sigma\left(  i\right)  }\right)  }_{=\sum_{P\subseteq\left\{
1,2,\ldots,n\right\}  }\left(  -1\right)  ^{\sigma}\left(  \prod_{i\in
P}a_{i,\sigma\left(  i\right)  }\right)  \left(  \prod_{i\in\widetilde{P}%
}b_{i,\sigma\left(  i\right)  }\right)  }\nonumber\\
&  =\underbrace{\sum_{\sigma\in S_{n}}\sum_{P\subseteq\left\{  1,2,\ldots
,n\right\}  }}_{=\sum_{P\subseteq\left\{  1,2,\ldots,n\right\}  }\sum
_{\sigma\in S_{n}}}\underbrace{\left(  -1\right)  ^{\sigma}\left(  \prod_{i\in
P}a_{i,\sigma\left(  i\right)  }\right)  \left(  \prod_{i\in\widetilde{P}%
}b_{i,\sigma\left(  i\right)  }\right)  }_{\substack{=c_{\sigma,P}\\\text{(by
(\ref{sol.det(A+B).c=}))}}}\nonumber\\
&  =\sum_{P\subseteq\left\{  1,2,\ldots,n\right\}  }\sum_{\sigma\in S_{n}%
}c_{\sigma,P}. \label{sol.det(A+B).1}%
\end{align}

But every subset $P$ of $\left\{  1,2,\ldots,n\right\}  $ satisfies%
\begin{align*}
&  \underbrace{\sum_{\sigma\in S_{n}}}_{\substack{=\sum_{Q\subseteq\left\{
1,2,\ldots,n\right\}  }\sum_{\substack{\sigma\in S_{n};\\\sigma\left(
P\right)  =Q}}\\\text{(because for every }\sigma\in S_{n}\text{, the
set}\\\sigma\left(  P\right)  \text{ is a subset of }\left\{  1,2,\ldots
,n\right\}  \text{)}}}c_{\sigma,P}\\
&  =\sum_{Q\subseteq\left\{  1,2,\ldots,n\right\}  }\sum_{\substack{\sigma\in
S_{n};\\\sigma\left(  P\right)  =Q}}c_{\sigma,P}\\
&  =\sum_{\substack{Q\subseteq\left\{  1,2,\ldots,n\right\}  ;\\\left\vert
Q\right\vert =\left\vert P\right\vert }}\sum_{\substack{\sigma\in
S_{n};\\\sigma\left(  P\right)  =Q}}c_{\sigma,P}+\sum_{\substack{Q\subseteq
\left\{  1,2,\ldots,n\right\}  ;\\\left\vert Q\right\vert \neq\left\vert
P\right\vert }}\underbrace{\sum_{\substack{\sigma\in S_{n};\\\sigma\left(
P\right)  =Q}}c_{\sigma,P}}_{\substack{=0\\\text{(by
(\ref{sol.det(A+B).emptysum}))}}}\\
&  \ \ \ \ \ \ \ \ \ \ \left(
\begin{array}
[c]{c}%
\text{since every subset }Q\text{ of }\left\{  1,2,\ldots,n\right\}  \text{
satisfies}\\
\text{either }\left\vert Q\right\vert =\left\vert P\right\vert \text{ or
}\left\vert Q\right\vert \neq\left\vert P\right\vert \text{ (but not both)}%
\end{array}
\right) \\
&  =\sum_{\substack{Q\subseteq\left\{  1,2,\ldots,n\right\}  ;\\\left\vert
Q\right\vert =\left\vert P\right\vert }}\sum_{\substack{\sigma\in
S_{n};\\\sigma\left(  P\right)  =Q}}c_{\sigma,P}+\underbrace{\sum
_{\substack{Q\subseteq\left\{  1,2,\ldots,n\right\}  ;\\\left\vert
Q\right\vert \neq\left\vert P\right\vert }}0}_{=0}\\
&  =\sum_{\substack{Q\subseteq\left\{  1,2,\ldots,n\right\}  ;\\\left\vert
Q\right\vert =\left\vert P\right\vert }}\underbrace{\sum_{\substack{\sigma\in
S_{n};\\\sigma\left(  P\right)  =Q}}c_{\sigma,P}}_{\substack{=\left(
-1\right)  ^{\sum P+\sum Q}\det\left(  \operatorname*{sub}\nolimits_{w\left(
P\right)  }^{w\left(  Q\right)  }A\right)  \det\left(  \operatorname*{sub}%
\nolimits_{w\left(  \widetilde{P}\right)  }^{w\left(  \widetilde{Q}\right)
}B\right)  \\\text{(by (\ref{sol.det(A+B).nonemptysum}))}}}\\
&  =\sum_{\substack{Q\subseteq\left\{  1,2,\ldots,n\right\}  ;\\\left\vert
Q\right\vert =\left\vert P\right\vert }}\left(  -1\right)  ^{\sum P+\sum
Q}\det\left(  \operatorname*{sub}\nolimits_{w\left(  P\right)  }^{w\left(
Q\right)  }A\right)  \det\left(  \operatorname*{sub}\nolimits_{w\left(
\widetilde{P}\right)  }^{w\left(  \widetilde{Q}\right)  }B\right)  .
\end{align*}
Hence, (\ref{sol.det(A+B).1}) becomes%
\begin{align*}
\det\left(  A+B\right)   &  =\sum_{P\subseteq\left\{  1,2,\ldots,n\right\}
}\underbrace{\sum_{\sigma\in S_{n}}c_{\sigma,P}}_{=\sum_{\substack{Q\subseteq
\left\{  1,2,\ldots,n\right\}  ;\\\left\vert Q\right\vert =\left\vert
P\right\vert }}\left(  -1\right)  ^{\sum P+\sum Q}\det\left(
\operatorname*{sub}\nolimits_{w\left(  P\right)  }^{w\left(  Q\right)
}A\right)  \det\left(  \operatorname*{sub}\nolimits_{w\left(  \widetilde{P}%
\right)  }^{w\left(  \widetilde{Q}\right)  }B\right)  }\\
&  =\sum_{P\subseteq\left\{  1,2,\ldots,n\right\}  }\sum_{\substack{Q\subseteq
\left\{  1,2,\ldots,n\right\}  ;\\\left\vert Q\right\vert =\left\vert
P\right\vert }}\left(  -1\right)  ^{\sum P+\sum Q}\det\left(
\operatorname*{sub}\nolimits_{w\left(  P\right)  }^{w\left(  Q\right)
}A\right)  \det\left(  \operatorname*{sub}\nolimits_{w\left(  \widetilde{P}%
\right)  }^{w\left(  \widetilde{Q}\right)  }B\right)  .
\end{align*}
This proves Theorem \ref{thm.det(A+B)}.
\end{proof}
\end{verlong}

\begin{proof}
[Solution to Exercise \ref{exe.det(A+B)}.]Exercise \ref{exe.det(A+B)} is
solved, since Theorem \ref{thm.det(A+B)} is proven.
\end{proof}

\subsection{Solution to Exercise \ref{exe.det(A+B).diag}}

\begin{proof}
[Proof of Lemma \ref{lem.diag.minors}.]The set $P$ is a subset of $\left\{
1,2,\ldots,n\right\}  $, and thus is finite (since $\left\{  1,2,\ldots
,n\right\}  $ is finite). Hence, $\left\vert P\right\vert \in\mathbb{N}$.

Define an element $k\in\mathbb{N}$ by $k=\left\vert P\right\vert $. Thus,
$k=\left\vert P\right\vert =\left\vert Q\right\vert $.

\begin{vershort}
The list $w\left(  P\right)  $ is the list of all elements of $P$ in
increasing order (with no repetitions), and thus has $k$ entries (since
$\left\vert P\right\vert =k$). Thus, write this list $w\left(  P\right)  $ in
the form $w\left(  P\right)  =\left(  p_{1},p_{2},\ldots,p_{k}\right)  $.
Similarly, write the list $w\left(  Q\right)  $ in the form $w\left(
Q\right)  =\left(  q_{1},q_{2},\ldots,q_{k}\right)  $.
\end{vershort}

\begin{verlong}
We know that $w\left(  P\right)  $ is the list of all elements of $P$ in
increasing order (with no repetitions) (by the definition of $w\left(
P\right)  $). Thus, $w\left(  P\right)  $ is a list of $\left\vert
P\right\vert $ elements. In other words, $w\left(  P\right)  $ is a list of
$k$ elements (since $\left\vert P\right\vert =k$).

Write $w\left(  P\right)  $ in the form $w\left(  P\right)  =\left(
p_{1},p_{2},\ldots,p_{k}\right)  $. (This is possible, since $w\left(
P\right)  $ is a list of $k$ elements.)

We know that $w\left(  Q\right)  $ is the list of all elements of $Q$ in
increasing order (with no repetitions) (by the definition of $w\left(
Q\right)  $). Thus, $w\left(  Q\right)  $ is a list of $\left\vert
Q\right\vert $ elements. In other words, $w\left(  Q\right)  $ is a list of
$k$ elements (since $\left\vert Q\right\vert =k$).

Write $w\left(  Q\right)  $ in the form $w\left(  Q\right)  =\left(
q_{1},q_{2},\ldots,q_{k}\right)  $. (This is possible, since $w\left(
Q\right)  $ is a list of $k$ elements.)
\end{verlong}

From $w\left(  P\right)  =\left(  p_{1},p_{2},\ldots,p_{k}\right)  $ and
$w\left(  Q\right)  =\left(  q_{1},q_{2},\ldots,q_{k}\right)  $, we obtain%
\begin{align}
\operatorname*{sub}\nolimits_{w\left(  P\right)  }^{w\left(  Q\right)  }D  &
=\operatorname*{sub}\nolimits_{\left(  p_{1},p_{2},\ldots,p_{k}\right)
}^{\left(  q_{1},q_{2},\ldots,q_{k}\right)  }D=\operatorname*{sub}%
\nolimits_{p_{1},p_{2},\ldots,p_{k}}^{q_{1},q_{2},\ldots,q_{k}}D\nonumber\\
&  =\left(  d_{p_{x}}\delta_{p_{x},q_{y}}\right)  _{1\leq x\leq k,\ 1\leq
y\leq k} \label{pf.lem.diag.minors.sub=}%
\end{align}
(by the definition of $\operatorname*{sub}\nolimits_{p_{1},p_{2},\ldots,p_{k}%
}^{q_{1},q_{2},\ldots,q_{k}}D$, since $D=\left(  d_{i}\delta_{i,j}\right)
_{1\leq i\leq n,\ 1\leq j\leq n}$). Thus, in particular, we see that
$\operatorname*{sub}\nolimits_{w\left(  P\right)  }^{w\left(  Q\right)  }D$ is
a $k\times k$-matrix.

\begin{vershort}
Recall that $\left(  p_{1},p_{2},\ldots,p_{k}\right)  =w\left(  P\right)  $ is
a list of all elements of $P$ (with no repetitions). Hence, $\left(
p_{1},p_{2},\ldots,p_{k}\right)  $ is a list with no repetitions, and
satisfies $\left\{  p_{1},p_{2},\ldots,p_{k}\right\}  =P$. Now, every $\left(
x,y\right)  \in\left\{  1,2,\ldots,k\right\}  ^{2}$ satisfies%
\begin{equation}
\delta_{p_{x},p_{y}}=\delta_{x,y} \label{pf.lem.diag.minors.short.deldel}%
\end{equation}
\footnote{\textit{Proof of (\ref{pf.lem.diag.minors.short.deldel}):} Let
$\left(  x,y\right)  \in\left\{  1,2,\ldots,k\right\}  ^{2}$. Then, the
integers $p_{1},p_{2},\ldots,p_{k}$ are pairwise distinct (since $\left(
p_{1},p_{2},\ldots,p_{k}\right)  $ is a list with no repetitions). Hence,
$p_{x}=p_{y}$ holds if and only if $x=y$. Thus, $%
\begin{cases}
1, & \text{if }p_{x}=p_{y};\\
0, & \text{if }p_{x}\neq p_{y}%
\end{cases}
=%
\begin{cases}
1, & \text{if }x=y;\\
0, & \text{if }x\neq y
\end{cases}
$. In view of $\delta_{p_{x},p_{y}}=%
\begin{cases}
1, & \text{if }p_{x}=p_{y};\\
0, & \text{if }p_{x}\neq p_{y}%
\end{cases}
$ and $\delta_{x,y}=%
\begin{cases}
1, & \text{if }x=y;\\
0, & \text{if }x\neq y
\end{cases}
$, this rewrites as $\delta_{p_{x},p_{y}}=\delta_{x,y}$. This proves
(\ref{pf.lem.diag.minors.short.deldel}).}.
\end{vershort}

\begin{verlong}
Notice that%
\begin{equation}
\delta_{p_{x},p_{y}}=\delta_{x,y}\ \ \ \ \ \ \ \ \ \ \text{for every }\left(
x,y\right)  \in\left\{  1,2,\ldots,k\right\}  ^{2}
\label{pf.lem.diag.minors.deldel}%
\end{equation}
\footnote{\textit{Proof of (\ref{pf.lem.diag.minors.deldel}):} Let $\left(
x,y\right)  \in\left\{  1,2,\ldots,k\right\}  ^{2}$. Thus, $x\in\left\{
1,2,\ldots,k\right\}  $ and $y\in\left\{  1,2,\ldots,k\right\}  $. We are in
one of the following two cases:
\par
\textit{Case 1:} We have $x=y$.
\par
\textit{Case 2:} We have $x\neq y$.
\par
Let us consider Case 1 first. In this case, we have $x=y$. Hence, $p_{x}%
=p_{y}$. Thus, $\delta_{p_{x},p_{y}}=1$. Comparing this with $\delta_{x,y}=1$
(since $x=y$), we obtain $\delta_{p_{x},p_{y}}=\delta_{x,y}$. Thus,
(\ref{pf.lem.diag.minors.deldel}) is proven in Case 1.
\par
Now, let us consider Case 2. In this case, we have $x\neq y$. In other words,
the elements $x$ and $y$ are distinct.
\par
But $w\left(  P\right)  $ is the list of all elements of $P$ in increasing
order (with no repetitions). Thus, $w\left(  P\right)  $ is a list with no
repetitions. In other words, $\left(  p_{1},p_{2},\ldots,p_{k}\right)  $ is a
list with no repetitions (since $w\left(  P\right)  =\left(  p_{1}%
,p_{2},\ldots,p_{k}\right)  $). In other words, the elements $p_{1}%
,p_{2},\ldots,p_{k}$ are pairwise distinct. In other words, if $u$ and $v$ are
any two distinct elements of $\left\{  1,2,\ldots,k\right\}  $, then
$p_{u}\neq p_{v}$. Applying this to $u=x$ and $v=y$, we obtain $p_{x}\neq
p_{y}$. Thus, $\delta_{p_{x},p_{y}}=0$. Comparing this with $\delta_{x,y}=0$
(since $x\neq y$), we obtain $\delta_{p_{x},p_{y}}=\delta_{x,y}$. Thus,
(\ref{pf.lem.diag.minors.deldel}) is proven in Case 2.
\par
We have now proven (\ref{pf.lem.diag.minors.deldel}) in each of the two Cases
1 and 2. Since these two Cases cover all possibilities, we thus conclude that
(\ref{pf.lem.diag.minors.deldel}) always holds.}.
\end{verlong}

\begin{vershort}
Recall again that $\left(  p_{1},p_{2},\ldots,p_{k}\right)  $ is a list of all
elements of $P$ (with no repetitions). Thus, the elements of $P$ are
$p_{1},p_{2},\ldots,p_{k}$ (with no repetitions). Hence,%
\begin{equation}
\prod_{i\in P}d_{i}=d_{p_{1}}d_{p_{2}}\cdots d_{p_{k}}=\prod_{i=1}^{k}%
d_{p_{i}}. \label{pf.lem.diag.minors.short.prod}%
\end{equation}

\end{vershort}

\begin{verlong}
Also,%
\begin{equation}
\prod_{i=1}^{k}d_{p_{i}}=\prod_{i\in P}d_{i} \label{pf.lem.diag.minors.prod}%
\end{equation}
\footnote{\textit{Proof of (\ref{pf.lem.diag.minors.prod}):} Let us use the
notation introduced in Definition \ref{def.sol.Ialbe.12n}. Thus, $\left[
k\right]  =\left\{  1,2,\ldots,k\right\}  $ (by the definition of $\left[
k\right]  $).
\par
We know that $w\left(  P\right)  $ is a list of all elements of $P$ (with no
repetitions). In other words, $\left(  p_{1},p_{2},\ldots,p_{k}\right)  $ is a
list of all elements of $P$ (with no repetitions) (since $w\left(  P\right)
=\left(  p_{1},p_{2},\ldots,p_{k}\right)  $). Thus, Lemma
\ref{lem.sol.Ialbe.inclist} \textbf{(a)} (applied to $P$, $k$ and $\left(
p_{1},p_{2},\ldots,p_{k}\right)  $ instead of $S$, $s$ and $\left(
c_{1},c_{2},\ldots,c_{s}\right)  $) shows that the map $\left[  k\right]
\rightarrow P,\ h\mapsto p_{h}$ is well-defined and a bijection. In
particular, this map is a bijection. Hence, we can substitute $p_{h}$ for $i$
in the product $\prod_{i\in P}d_{i}$. We thus obtain%
\[
\prod_{i\in P}d_{i}=\underbrace{\prod_{h\in\left[  k\right]  }}%
_{\substack{=\prod_{h\in\left\{  1,2,\ldots,k\right\}  }\\\text{(since
}\left[  k\right]  =\left\{  1,2,\ldots,k\right\}  \text{)}}}d_{p_{h}%
}=\underbrace{\prod_{h\in\left\{  1,2,\ldots,k\right\}  }}_{=\prod_{h=1}^{k}%
}d_{p_{h}}=\prod_{h=1}^{k}d_{p_{h}}=\prod_{i=1}^{k}d_{p_{i}}%
\]
(here, we have renamed the index $h$ as $i$ in the product). This proves
(\ref{pf.lem.diag.minors.prod}).}.
\end{verlong}

We are in one of the following two cases:

\textit{Case 1:} We have $P\neq Q$.

\textit{Case 2:} We have $P=Q$.

\begin{vershort}
Let us consider Case 1 first. In this case, we have $P\neq Q$. Hence,
$\delta_{P,Q}=0$. On the other hand, there exists some $p\in P$ such that
$p\notin Q$\ \ \ \ \footnote{\textit{Proof.} Assume the contrary (for the sake
of contradiction). Hence, every $p\in P$ satisfies $p\in Q$. In other words,
$P\subseteq Q$. But $\left\vert P\right\vert =\left\vert Q\right\vert
\geq\left\vert Q\right\vert $.
\par
But if $X$ is a finite set, then every subset of $X$ that has size
$\geq\left\vert X\right\vert $ must be $X$ itself. In other words, if $X$ is a
finite set, and if $Y$ is a subset of $X$ satisfying $\left\vert Y\right\vert
\geq\left\vert X\right\vert $, then $Y=X$. Applying this to $X=Q$ and $Y=P$,
we conclude that $P=Q$ (since $P$ is a subset of $Q$, and since $Q$ is a
finite set). This contradicts $P\neq Q$. This contradiction shows that our
assumption was wrong; qed.}. Consider this $p$.
\end{vershort}

\begin{verlong}
Let us consider Case 1 first. In this case, we have $P\neq Q$. Hence,
$\delta_{P,Q}=0$. On the other hand, there exists some $p\in P$ such that
$p\notin Q$\ \ \ \ \footnote{\textit{Proof.} Assume the contrary (for the sake
of contradiction). Hence, every $p\in P$ satisfies $p\in Q$. In other words,
$P\subseteq Q$. In other words, $P$ is a subset of $Q$. But $\left\vert
P\right\vert =\left\vert Q\right\vert \geq\left\vert Q\right\vert $.
\par
We know that $Q$ is a subset of $\left\{  1,2,\ldots,n\right\}  $, and thus is
a finite set (since $\left\{  1,2,\ldots,n\right\}  $ is a finite set).
\par
But if $X$ is a finite set, then every subset of $X$ that has size
$\geq\left\vert X\right\vert $ must be $X$ itself. In other words, if $X$ is a
finite set, and if $Y$ is a subset of $X$ satisfying $\left\vert Y\right\vert
\geq\left\vert X\right\vert $, then $Y=X$. Applying this to $X=Q$ and $Y=P$,
we conclude that $P=Q$ (since $P$ is a subset of $Q$). This contradicts $P\neq
Q$. This contradiction shows that our assumption was wrong; qed.}. Consider
this $p$.
\end{verlong}

\begin{vershort}
Now, $p\in P=\left\{  p_{1},p_{2},\ldots,p_{k}\right\}  $. In other words,
$p=p_{u}$ for some $u\in\left\{  1,2,\ldots,k\right\}  $. Consider this $u$.
Every $y\in\left\{  1,2,\ldots,k\right\}  $ satisfies%
\begin{equation}
\delta_{p_{u},q_{y}}=0 \label{pf.lem.diag.minors.short.C1.1}%
\end{equation}
\footnote{\textit{Proof of (\ref{pf.lem.diag.minors.short.C1.1}):} Let
$y\in\left\{  1,2,\ldots,k\right\}  $.
\par
We know that $\left(  q_{1},q_{2},\ldots,q_{k}\right)  =w\left(  Q\right)  $
is a list of all elements of $Q$. Hence, $\left\{  q_{1},q_{2},\ldots
,q_{k}\right\}  =Q$. Now, $y\in\left\{  1,2,\ldots,k\right\}  $, so that
$q_{y}\in\left\{  q_{1},q_{2},\ldots,q_{k}\right\}  =Q$. If we had
$p_{u}=q_{y}$, then we would have $p=p_{u}=q_{y}\in Q$, which would contradict
$p\notin Q$. Thus, we cannot have $p_{u}=q_{y}$. In other words, we have
$p_{u}\neq q_{y}$. Hence, $\delta_{p_{u},q_{y}}=0$. This proves
(\ref{pf.lem.diag.minors.short.C1.1}).}. Now,
\begin{align*}
&  \left(  \text{the }u\text{-th row of the matrix }%
\underbrace{\operatorname*{sub}\nolimits_{w\left(  P\right)  }^{w\left(
Q\right)  }D}_{=\left(  d_{p_{x}}\delta_{p_{x},q_{y}}\right)  _{1\leq x\leq
k,\ 1\leq y\leq k}}\right) \\
&  =\left(  \text{the }u\text{-th row of the matrix }\left(  d_{p_{x}}%
\delta_{p_{x},q_{y}}\right)  _{1\leq x\leq k,\ 1\leq y\leq k}\right) \\
&  =\left(  d_{p_{u}}\underbrace{\delta_{p_{u},q_{y}}}%
_{\substack{=0\\\text{(by (\ref{pf.lem.diag.minors.short.C1.1}))}}}\right)
_{1\leq x\leq1,\ 1\leq y\leq k}=\left(  \underbrace{d_{p_{u}}0}_{=0}\right)
_{1\leq x\leq1,\ 1\leq y\leq k}=\left(  0\right)  _{1\leq x\leq1,\ 1\leq y\leq
k}.
\end{align*}
In other words, the $u$-th row of the matrix $\operatorname*{sub}%
\nolimits_{w\left(  P\right)  }^{w\left(  Q\right)  }D$ consists of zeroes.
Therefore, Exercise \ref{exe.ps4.6} \textbf{(c)} (applied to $k$ and
$\operatorname*{sub}\nolimits_{w\left(  P\right)  }^{w\left(  Q\right)  }D$
instead of $n$ and $A$) yields $\det\left(  \operatorname*{sub}%
\nolimits_{w\left(  P\right)  }^{w\left(  Q\right)  }D\right)  =0$. Comparing
this with $\underbrace{\delta_{P,Q}}_{=0}\prod_{i\in P}d_{i}=0$, we obtain
$\det\left(  \operatorname*{sub}\nolimits_{w\left(  P\right)  }^{w\left(
Q\right)  }D\right)  =\delta_{P,Q}\prod_{i\in P}d_{i}$. Hence, Lemma
\ref{lem.diag.minors} is proven in Case 1.
\end{vershort}

\begin{verlong}
We know that $w\left(  P\right)  $ is a list of all elements of $P$. In other
words, $\left(  p_{1},p_{2},\ldots,p_{k}\right)  $ is a list of all elements
of $P$ (since $w\left(  P\right)  =\left(  p_{1},p_{2},\ldots,p_{k}\right)
$). Hence, $\left\{  p_{1},p_{2},\ldots,p_{k}\right\}  =P$. Now, $p\in
P=\left\{  p_{1},p_{2},\ldots,p_{k}\right\}  $. In other words, $p=p_{u}$ for
some $u\in\left\{  1,2,\ldots,k\right\}  $. Consider this $u$. Every
$y\in\left\{  1,2,\ldots,k\right\}  $ satisfies%
\begin{equation}
\delta_{p_{u},q_{y}}=0 \label{pf.lem.diag.minors.C1.1}%
\end{equation}
\footnote{\textit{Proof of (\ref{pf.lem.diag.minors.C1.1}):} Let $y\in\left\{
1,2,\ldots,k\right\}  $.
\par
We know that $w\left(  Q\right)  $ is a list of all elements of $Q$. In other
words, $\left(  q_{1},q_{2},\ldots,q_{k}\right)  $ is a list of all elements
of $Q$ (since $w\left(  Q\right)  =\left(  q_{1},q_{2},\ldots,q_{k}\right)
$). Hence, $\left\{  q_{1},q_{2},\ldots,q_{k}\right\}  =Q$. Now, $y\in\left\{
1,2,\ldots,k\right\}  $, so that $q_{y}\in\left\{  q_{1},q_{2},\ldots
,q_{k}\right\}  =Q$. If we had $p_{u}=q_{y}$, then we would have
$p=p_{u}=q_{y}\in Q$, which would contradict $p\notin Q$. Thus, we cannot have
$p_{u}=q_{y}$. In other words, we have $p_{u}\neq q_{y}$. Hence,
$\delta_{p_{u},q_{y}}=0$. This proves (\ref{pf.lem.diag.minors.C1.1}).}. Now,
\begin{align*}
&  \left(  \text{the }u\text{-th row of the matrix }%
\underbrace{\operatorname*{sub}\nolimits_{w\left(  P\right)  }^{w\left(
Q\right)  }D}_{\substack{=\left(  d_{p_{x}}\delta_{p_{x},q_{y}}\right)
_{1\leq x\leq k,\ 1\leq y\leq k}\\\text{(by (\ref{pf.lem.diag.minors.sub=}))}%
}}\right) \\
&  =\left(  \text{the }u\text{-th row of the matrix }\left(  d_{p_{x}}%
\delta_{p_{x},q_{y}}\right)  _{1\leq x\leq k,\ 1\leq y\leq k}\right) \\
&  =\left(  d_{p_{u}}\underbrace{\delta_{p_{u},q_{y}}}%
_{\substack{=0\\\text{(by (\ref{pf.lem.diag.minors.C1.1}))}}}\right)  _{1\leq
x\leq1,\ 1\leq y\leq k}=\left(  \underbrace{d_{p_{u}}0}_{=0}\right)  _{1\leq
x\leq1,\ 1\leq y\leq k}=\left(  0\right)  _{1\leq x\leq1,\ 1\leq y\leq k}.
\end{align*}
In other words, the $u$-th row of the matrix $\operatorname*{sub}%
\nolimits_{w\left(  P\right)  }^{w\left(  Q\right)  }D$ consists of zeroes.
Hence, a row of the matrix $\operatorname*{sub}\nolimits_{w\left(  P\right)
}^{w\left(  Q\right)  }D$ consists of zeroes (namely, the $u$-th row). Hence,
Exercise \ref{exe.ps4.6} \textbf{(c)} (applied to $k$ and $\operatorname*{sub}%
\nolimits_{w\left(  P\right)  }^{w\left(  Q\right)  }D$ instead of $n$ and
$A$) yields $\det\left(  \operatorname*{sub}\nolimits_{w\left(  P\right)
}^{w\left(  Q\right)  }D\right)  =0$. Comparing this with $\underbrace{\delta
_{P,Q}}_{=0}\prod_{i\in P}d_{i}=0$, we obtain $\det\left(  \operatorname*{sub}%
\nolimits_{w\left(  P\right)  }^{w\left(  Q\right)  }D\right)  =\delta
_{P,Q}\prod_{i\in P}d_{i}$. Hence, Lemma \ref{lem.diag.minors} is proven in
Case 1.
\end{verlong}

\begin{vershort}
Let us now consider Case 2. In this case, we have $P=Q$. Hence, $\delta
_{P,Q}=1$. On the other hand, from $w\left(  Q\right)  =\left(  q_{1}%
,q_{2},\ldots,q_{k}\right)  $, we obtain
\[
\left(  q_{1},q_{2},\ldots,q_{k}\right)  =w\left(  \underbrace{Q}_{=P}\right)
=w\left(  P\right)  =\left(  p_{1},p_{2},\ldots,p_{k}\right)  .
\]
In other words, $q_{y}=p_{y}$ for every $y\in\left\{  1,2,\ldots,k\right\}  $.
Thus, every $y\in\left\{  1,2,\ldots,k\right\}  $ satisfies%
\begin{align}
\delta_{p_{x},q_{y}}  &  =\delta_{p_{x},p_{y}}\ \ \ \ \ \ \ \ \ \ \left(
\text{since }q_{y}=p_{y}\right) \nonumber\\
&  =\delta_{x,y}\ \ \ \ \ \ \ \ \ \ \left(  \text{by
(\ref{pf.lem.diag.minors.short.deldel})}\right)  .
\label{pf.lem.diag.minors.short.C2.2}%
\end{align}
Now, (\ref{pf.lem.diag.minors.sub=}) becomes%
\[
\operatorname*{sub}\nolimits_{w\left(  P\right)  }^{w\left(  Q\right)
}D=\left(  d_{p_{x}}\underbrace{\delta_{p_{x},q_{y}}}_{\substack{=\delta
_{x,y}\\\text{(by (\ref{pf.lem.diag.minors.short.C2.2}))}}}\right)  _{1\leq
x\leq k,\ 1\leq y\leq k}=\left(  d_{p_{x}}\delta_{x,y}\right)  _{1\leq x\leq
k,\ 1\leq y\leq k}=\left(  d_{p_{i}}\delta_{i,j}\right)  _{1\leq i\leq
k,\ 1\leq j\leq k}%
\]
(here, we have renamed the index $\left(  x,y\right)  $ as $\left(
i,j\right)  $).
\end{vershort}

\begin{verlong}
Let us now consider Case 2. In this case, we have $P=Q$. Hence, $\delta
_{P,Q}=1$. On the other hand, from $w\left(  Q\right)  =\left(  q_{1}%
,q_{2},\ldots,q_{k}\right)  $, we obtain
\[
\left(  q_{1},q_{2},\ldots,q_{k}\right)  =w\left(  \underbrace{Q}_{=P}\right)
=w\left(  P\right)  =\left(  p_{1},p_{2},\ldots,p_{k}\right)  .
\]
In other words,%
\begin{equation}
q_{y}=p_{y}\ \ \ \ \ \ \ \ \ \ \text{for every }y\in\left\{  1,2,\ldots
,k\right\}  . \label{pf.lem.diag.minors.C2.1}%
\end{equation}
Now, (\ref{pf.lem.diag.minors.sub=}) becomes%
\begin{align*}
\operatorname*{sub}\nolimits_{w\left(  P\right)  }^{w\left(  Q\right)  }D  &
=\left(  d_{p_{x}}\underbrace{\delta_{p_{x},q_{y}}}_{\substack{=\delta
_{p_{x},p_{y}}\\\text{(by (\ref{pf.lem.diag.minors.C2.1}))}}}\right)  _{1\leq
x\leq k,\ 1\leq y\leq k}=\left(  d_{p_{x}}\underbrace{\delta_{p_{x},p_{y}}%
}_{\substack{=\delta_{x,y}\\\text{(by (\ref{pf.lem.diag.minors.deldel}))}%
}}\right)  _{1\leq x\leq k,\ 1\leq y\leq k}\\
&  =\left(  d_{p_{x}}\delta_{x,y}\right)  _{1\leq x\leq k,\ 1\leq y\leq
k}=\left(  d_{p_{i}}\delta_{i,j}\right)  _{1\leq i\leq k,\ 1\leq j\leq k}%
\end{align*}
(here, we have renamed the index $\left(  x,y\right)  $ as $\left(
i,j\right)  $).
\end{verlong}

\begin{vershort}
But we have $d_{p_{i}}\delta_{i,j}=0$ for every $\left(  i,j\right)
\in\left\{  1,2,\ldots,k\right\}  ^{2}$ satisfying $i<j$%
\ \ \ \ \footnote{\textit{Proof.} Let $\left(  i,j\right)  \in\left\{
1,2,\ldots,k\right\}  ^{2}$ be such that $i<j$. Thus, $i\neq j$ (since $i<j$),
so that $\delta_{i,j}=0$. Thus, $d_{p_{i}}\underbrace{\delta_{i,j}}_{=0}=0$,
qed.}. Hence, Exercise \ref{exe.ps4.3} (applied to $k$, $\operatorname*{sub}%
\nolimits_{w\left(  P\right)  }^{w\left(  Q\right)  }D$ and $d_{p_{i}}%
\delta_{i,j}$ instead of $n$, $A$ and $a_{i,j}$) yields $\det\left(
\operatorname*{sub}\nolimits_{w\left(  P\right)  }^{w\left(  Q\right)
}D\right)  =\left(  d_{p_{1}}\delta_{1,1}\right)  \left(  d_{p_{2}}%
\delta_{2,2}\right)  \cdots\left(  d_{p_{k}}\delta_{k,k}\right)  $ (since
$\operatorname*{sub}\nolimits_{w\left(  P\right)  }^{w\left(  Q\right)
}D=\left(  d_{p_{i}}\delta_{i,j}\right)  _{1\leq i\leq k,\ 1\leq j\leq k}$).
Thus,%
\begin{align*}
\det\left(  \operatorname*{sub}\nolimits_{w\left(  P\right)  }^{w\left(
Q\right)  }D\right)   &  =\left(  d_{p_{1}}\delta_{1,1}\right)  \left(
d_{p_{2}}\delta_{2,2}\right)  \cdots\left(  d_{p_{k}}\delta_{k,k}\right)
=\prod_{i=1}^{k}\left(  d_{p_{i}}\underbrace{\delta_{i,i}}%
_{\substack{=1\\\text{(since }i=i\text{)}}}\right) \\
&  =\prod_{i=1}^{k}d_{p_{i}}=\prod_{i\in P}d_{i}\ \ \ \ \ \ \ \ \ \ \left(
\text{by (\ref{pf.lem.diag.minors.short.prod})}\right)  .
\end{align*}
Comparing this with $\underbrace{\delta_{P,Q}}_{=1}\prod_{i\in P}d_{i}%
=\prod_{i\in P}d_{i}$, we obtain $\det\left(  \operatorname*{sub}%
\nolimits_{w\left(  P\right)  }^{w\left(  Q\right)  }D\right)  =\delta
_{P,Q}\prod_{i\in P}d_{i}$. Hence, Lemma \ref{lem.diag.minors} is proven in
Case 2.
\end{vershort}

\begin{verlong}
But we have $d_{p_{i}}\delta_{i,j}=0$ for every $\left(  i,j\right)
\in\left\{  1,2,\ldots,k\right\}  ^{2}$ satisfying $i<j$%
\ \ \ \ \footnote{\textit{Proof.} Let $\left(  i,j\right)  \in\left\{
1,2,\ldots,k\right\}  ^{2}$ be such that $i<j$. Thus, $i\neq j$ (since $i<j$),
so that $\delta_{i,j}=0$. Thus, $d_{p_{i}}\underbrace{\delta_{i,j}}_{=0}=0$,
qed.}. Hence, Exercise \ref{exe.ps4.3} (applied to $k$, $\operatorname*{sub}%
\nolimits_{w\left(  P\right)  }^{w\left(  Q\right)  }D$ and $d_{p_{i}}%
\delta_{i,j}$ instead of $n$, $A$ and $a_{i,j}$) yields $\det\left(
\operatorname*{sub}\nolimits_{w\left(  P\right)  }^{w\left(  Q\right)
}D\right)  =\left(  d_{p_{1}}\delta_{1,1}\right)  \left(  d_{p_{2}}%
\delta_{2,2}\right)  \cdots\left(  d_{p_{k}}\delta_{k,k}\right)  $ (since
$\operatorname*{sub}\nolimits_{w\left(  P\right)  }^{w\left(  Q\right)
}D=\left(  d_{p_{i}}\delta_{i,j}\right)  _{1\leq i\leq k,\ 1\leq j\leq k}$).
Thus,%
\begin{align*}
\det\left(  \operatorname*{sub}\nolimits_{w\left(  P\right)  }^{w\left(
Q\right)  }D\right)   &  =\left(  d_{p_{1}}\delta_{1,1}\right)  \left(
d_{p_{2}}\delta_{2,2}\right)  \cdots\left(  d_{p_{k}}\delta_{k,k}\right)
=\prod_{i=1}^{k}\left(  d_{p_{i}}\underbrace{\delta_{i,i}}%
_{\substack{=1\\\text{(since }i=i\text{)}}}\right) \\
&  =\prod_{i=1}^{k}d_{p_{i}}=\prod_{i\in P}d_{i}\ \ \ \ \ \ \ \ \ \ \left(
\text{by (\ref{pf.lem.diag.minors.prod})}\right)  .
\end{align*}
Comparing this with $\underbrace{\delta_{P,Q}}_{=1}\prod_{i\in P}d_{i}%
=\prod_{i\in P}d_{i}$, we obtain $\det\left(  \operatorname*{sub}%
\nolimits_{w\left(  P\right)  }^{w\left(  Q\right)  }D\right)  =\delta
_{P,Q}\prod_{i\in P}d_{i}$. Hence, Lemma \ref{lem.diag.minors} is proven in
Case 2.
\end{verlong}

We now have proven Lemma \ref{lem.diag.minors} in each of the two Cases 1 and
2. Since these two Cases cover all possibilities, we thus conclude that Lemma
\ref{lem.diag.minors} always holds.
\end{proof}

\begin{proof}
[Proof of Corollary \ref{cor.det(A+D)}.]We start by making some general observations:

\begin{vershort}
For any subset $I$ of $\left\{  1,2,\ldots,n\right\}  $, we let $\widetilde{I}%
$ denote the complement $\left\{  1,2,\ldots,n\right\}  \setminus I$ of $I$.
Then, any two subsets $P$ and $Q$ of $\left\{  1,2,\ldots,n\right\}  $ satisfy%
\begin{equation}
\delta_{\widetilde{P},\widetilde{Q}}=\delta_{P,Q}
\label{pf.cor.det(A+D).short.del=del}%
\end{equation}
\footnote{\textit{Proof of (\ref{pf.cor.det(A+D).short.del=del}):} Let $P$ and
$Q$ be two subsets of $\left\{  1,2,\ldots,n\right\}  $. If $P=Q$, then
$\widetilde{P}=\widetilde{Q}$ and thus $\delta_{\widetilde{P},\widetilde{Q}%
}=1=\delta_{P,Q}$ (since $P=Q$ yields $\delta_{P,Q}=1$). Therefore, if $P=Q$,
then (\ref{pf.cor.det(A+D).short.del=del}) holds. Hence, for the rest of the
proof of (\ref{pf.cor.det(A+D).short.del=del}), we WLOG assume that we don't
have $P=Q$. In other words, we have $P\neq Q$. Thus, $\delta_{P,Q}=0$.
\par
Now, the definition of $\widetilde{P}$ yields $\widetilde{P}=\left\{
1,2,\ldots,n\right\}  \setminus P$. Hence, $\widetilde{P}$ is again a subset
of $\left\{  1,2,\ldots,n\right\}  $. The definition of
$\widetilde{\widetilde{P}}$ yields
\[
\widetilde{\widetilde{P}}=\left\{  1,2,\ldots,n\right\}  \setminus
\underbrace{\widetilde{P}}_{=\left\{  1,2,\ldots,n\right\}  \setminus
P}=\left\{  1,2,\ldots,n\right\}  \setminus\left(  \left\{  1,2,\ldots
,n\right\}  \setminus P\right)  =P
\]
(since $P\subseteq\left\{  1,2,\ldots,n\right\}  $). Similarly,
$\widetilde{\widetilde{Q}}=Q$. If we had $\widetilde{P}=\widetilde{Q}$, then
we would have $\widetilde{\widetilde{P}}=\widetilde{\widetilde{Q}}$, which
would contradict $\widetilde{\widetilde{P}}=P\neq Q=\widetilde{\widetilde{Q}}%
$. Hence, we cannot have $\widetilde{P}=\widetilde{Q}$. Thus, we have
$\widetilde{P}\neq\widetilde{Q}$, so that $\delta_{\widetilde{P}%
,\widetilde{Q}}=0$. Compared with $\delta_{P,Q}=0$, this yields $\delta
_{\widetilde{P},\widetilde{Q}}=\delta_{P,Q}$. This proves
(\ref{pf.cor.det(A+D).short.del=del}).}.
\end{vershort}

\begin{verlong}
For any subset $I$ of $\left\{  1,2,\ldots,n\right\}  $, we let $\widetilde{I}%
$ denote the complement $\left\{  1,2,\ldots,n\right\}  \setminus I$ of $I$.
Then, any two subsets $P$ and $Q$ of $\left\{  1,2,\ldots,n\right\}  $ satisfy%
\begin{equation}
\delta_{\widetilde{P},\widetilde{Q}}=\delta_{P,Q}
\label{pf.cor.det(A+D).del=del}%
\end{equation}
\footnote{\textit{Proof of (\ref{pf.cor.det(A+D).del=del}):} Let $P$ and $Q$
be two subsets of $\left\{  1,2,\ldots,n\right\}  $.
\par
The definition of $\widetilde{P}$ yields $\widetilde{P}=\left\{
1,2,\ldots,n\right\}  \setminus P\subseteq\left\{  1,2,\ldots,n\right\}  $.
Hence, $\widetilde{P}$ is again a subset of $\left\{  1,2,\ldots,n\right\}  $.
Thus, the definition of $\widetilde{\widetilde{P}}$ yields
\[
\widetilde{\widetilde{P}}=\left\{  1,2,\ldots,n\right\}  \setminus
\underbrace{\widetilde{P}}_{=\left\{  1,2,\ldots,n\right\}  \setminus
P}=\left\{  1,2,\ldots,n\right\}  \setminus\left(  \left\{  1,2,\ldots
,n\right\}  \setminus P\right)  =P
\]
(since $P\subseteq\left\{  1,2,\ldots,n\right\}  $). The same argument
(applied to $Q$ instead of $P$) shows that $\widetilde{\widetilde{Q}}=Q$.
\par
If $P=Q$, then%
\begin{align*}
\delta_{\widetilde{P},\widetilde{Q}}  &  =\delta_{\widetilde{Q},\widetilde{Q}%
}\ \ \ \ \ \ \ \ \ \ \left(  \text{since }P=Q\right) \\
&  =1\ \ \ \ \ \ \ \ \ \ \left(  \text{since }\widetilde{Q}=\widetilde{Q}%
\right) \\
&  =\delta_{P,Q}\ \ \ \ \ \ \ \ \ \ \left(  \text{since }\delta_{P,Q}=1\text{
(since }P=Q\text{)}\right)  .
\end{align*}
Hence, if $P=Q$, then (\ref{pf.cor.det(A+D).del=del}) is proven. Thus, for the
rest of our proof of (\ref{pf.cor.det(A+D).del=del}), we can WLOG assume that
we don't have $P=Q$. Assume this.
\par
We don't have $P=Q$. In other words, we have $P\neq Q$. Hence, $\delta
_{P,Q}=0$. If we had $\widetilde{P}=\widetilde{Q}$, then we would have
\begin{align*}
\widetilde{\widetilde{P}}  &  =\widetilde{\widetilde{Q}}%
\ \ \ \ \ \ \ \ \ \ \left(  \text{since }\widetilde{P}=\widetilde{Q}\right) \\
&  =Q,
\end{align*}
which would contradict $\widetilde{\widetilde{P}}=P\neq Q$. Hence, we cannot
have $\widetilde{P}=\widetilde{Q}$. In other words, we have $\widetilde{P}%
\neq\widetilde{Q}$. Hence, $\delta_{\widetilde{P},\widetilde{Q}}=0$. Comparing
this with $\delta_{P,Q}=0$ (since $P\neq Q$), we obtain $\delta_{\widetilde{P}%
,\widetilde{Q}}=\delta_{P,Q}$. This proves (\ref{pf.cor.det(A+D).del=del}).}.
\end{verlong}

\begin{vershort}
Furthermore, if $P$ and $Q$ are two subsets of $\left\{  1,2,\ldots,n\right\}
$ satisfying $\left\vert P\right\vert =\left\vert Q\right\vert $, then%
\begin{equation}
\det\left(  \operatorname*{sub}\nolimits_{w\left(  \widetilde{P}\right)
}^{w\left(  \widetilde{Q}\right)  }D\right)  =\delta_{P,Q}\prod_{i\in\left\{
1,2,\ldots,n\right\}  \setminus P}d_{i} \label{pf.cor.det(A+D).short.1}%
\end{equation}
\footnote{\textit{Proof of (\ref{pf.cor.det(A+D).short.1}):} Let $P$ and $Q$
be two subsets of $\left\{  1,2,\ldots,n\right\}  $ satisfying $\left\vert
P\right\vert =\left\vert Q\right\vert $.
\par
The definition of $\widetilde{P}$ yields $\widetilde{P}=\left\{
1,2,\ldots,n\right\}  \setminus P\subseteq\left\{  1,2,\ldots,n\right\}  $. In
other words, $\widetilde{P}$ is a subset of $\left\{  1,2,\ldots,n\right\}  $.
The same argument (applied to $Q$ instead of $P$) shows that $\widetilde{Q}$
is a subset of $\left\{  1,2,\ldots,n\right\}  $.
\par
We have%
\begin{align*}
\left\vert \underbrace{\widetilde{P}}_{=\left\{  1,2,\ldots,n\right\}
\setminus P}\right\vert  &  =\left\vert \left\{  1,2,\ldots,n\right\}
\setminus P\right\vert =\underbrace{\left\vert \left\{  1,2,\ldots,n\right\}
\right\vert }_{=n}-\left\vert P\right\vert \ \ \ \ \ \ \ \ \ \ \left(
\text{since }P\subseteq\left\{  1,2,\ldots,n\right\}  \right) \\
&  =n-\left\vert P\right\vert .
\end{align*}
The same argument (applied to $Q$ instead of $P$) shows that $\left\vert
\widetilde{Q}\right\vert =n-\left\vert Q\right\vert $. Now, $\left\vert
\widetilde{P}\right\vert =n-\underbrace{\left\vert P\right\vert }_{=\left\vert
Q\right\vert }=n-\left\vert Q\right\vert =\left\vert \widetilde{Q}\right\vert
$. Thus, Lemma \ref{lem.diag.minors} (applied to $\widetilde{P}$ and
$\widetilde{Q}$ instead of $P$ and $Q$) yields
\[
\det\left(  \operatorname*{sub}\nolimits_{w\left(  \widetilde{P}\right)
}^{w\left(  \widetilde{Q}\right)  }D\right)  =\underbrace{\delta
_{\widetilde{P},\widetilde{Q}}}_{\substack{=\delta_{P,Q}\\\text{(by
(\ref{pf.cor.det(A+D).short.del=del}))}}}\underbrace{\prod_{i\in\widetilde{P}%
}}_{\substack{=\prod_{i\in\left\{  1,2,\ldots,n\right\}  \setminus
P}\\\text{(since }\widetilde{P}=\left\{  1,2,\ldots,n\right\}  \setminus
P\text{)}}}d_{i}=\delta_{P,Q}\prod_{i\in\left\{  1,2,\ldots,n\right\}
\setminus P}d_{i}.
\]
This proves (\ref{pf.cor.det(A+D).short.1}).}.
\end{vershort}

\begin{verlong}
Furthermore, if $P$ and $Q$ are two subsets of $\left\{  1,2,\ldots,n\right\}
$ satisfying $\left\vert P\right\vert =\left\vert Q\right\vert $, then%
\begin{equation}
\det\left(  \operatorname*{sub}\nolimits_{w\left(  \widetilde{P}\right)
}^{w\left(  \widetilde{Q}\right)  }D\right)  =\delta_{P,Q}\prod_{i\in\left\{
1,2,\ldots,n\right\}  \setminus P}d_{i} \label{pf.cor.det(A+D).1}%
\end{equation}
\footnote{\textit{Proof of (\ref{pf.cor.det(A+D).1}):} Let $P$ and $Q$ be two
subsets of $\left\{  1,2,\ldots,n\right\}  $ satisfying $\left\vert
P\right\vert =\left\vert Q\right\vert $.
\par
The definition of $\widetilde{P}$ yields $\widetilde{P}=\left\{
1,2,\ldots,n\right\}  \setminus P\subseteq\left\{  1,2,\ldots,n\right\}  $. In
other words, $\widetilde{P}$ is a subset of $\left\{  1,2,\ldots,n\right\}  $.
The same argument (applied to $Q$ instead of $P$) shows that $\widetilde{Q}$
is a subset of $\left\{  1,2,\ldots,n\right\}  $.
\par
We have%
\begin{align*}
\left\vert \underbrace{\widetilde{P}}_{=\left\{  1,2,\ldots,n\right\}
\setminus P}\right\vert  &  =\left\vert \left\{  1,2,\ldots,n\right\}
\setminus P\right\vert =\underbrace{\left\vert \left\{  1,2,\ldots,n\right\}
\right\vert }_{=n}-\left\vert P\right\vert \ \ \ \ \ \ \ \ \ \ \left(
\text{since }P\subseteq\left\{  1,2,\ldots,n\right\}  \right) \\
&  =n-\left\vert P\right\vert .
\end{align*}
The same argument (applied to $Q$ instead of $P$) shows that $\left\vert
\widetilde{Q}\right\vert =n-\left\vert Q\right\vert $. Now, $\left\vert
\widetilde{P}\right\vert =n-\underbrace{\left\vert P\right\vert }_{=\left\vert
Q\right\vert }=n-\left\vert Q\right\vert =\left\vert \widetilde{Q}\right\vert
$. Thus, Lemma \ref{lem.diag.minors} (applied to $\widetilde{P}$ and
$\widetilde{Q}$ instead of $P$ and $Q$) yields
\[
\det\left(  \operatorname*{sub}\nolimits_{w\left(  \widetilde{P}\right)
}^{w\left(  \widetilde{Q}\right)  }D\right)  =\underbrace{\delta
_{\widetilde{P},\widetilde{Q}}}_{\substack{=\delta_{P,Q}\\\text{(by
(\ref{pf.cor.det(A+D).del=del}))}}}\underbrace{\prod_{i\in\widetilde{P}}%
}_{\substack{=\prod_{i\in\left\{  1,2,\ldots,n\right\}  \setminus
P}\\\text{(since }\widetilde{P}=\left\{  1,2,\ldots,n\right\}  \setminus
P\text{)}}}d_{i}=\delta_{P,Q}\prod_{i\in\left\{  1,2,\ldots,n\right\}
\setminus P}d_{i}.
\]
This proves (\ref{pf.cor.det(A+D).1}).}.
\end{verlong}

\begin{vershort}
Now, every subset $P$ of $\left\{  1,2,\ldots,n\right\}  $ satisfies%
\begin{align}
&  \sum_{\substack{Q\subseteq\left\{  1,2,\ldots,n\right\}  ;\\\left\vert
P\right\vert =\left\vert Q\right\vert }}\left(  -1\right)  ^{\sum P+\sum
Q}\det\left(  \operatorname*{sub}\nolimits_{w\left(  P\right)  }^{w\left(
Q\right)  }A\right)  \underbrace{\det\left(  \operatorname*{sub}%
\nolimits_{w\left(  \widetilde{P}\right)  }^{w\left(  \widetilde{Q}\right)
}D\right)  }_{\substack{=\delta_{P,Q}\prod_{i\in\left\{  1,2,\ldots,n\right\}
\setminus P}d_{i}\\\text{(by (\ref{pf.cor.det(A+D).short.1}))}}}\nonumber\\
&  =\sum_{\substack{Q\subseteq\left\{  1,2,\ldots,n\right\}  ;\\\left\vert
P\right\vert =\left\vert Q\right\vert }}\left(  -1\right)  ^{\sum P+\sum
Q}\det\left(  \operatorname*{sub}\nolimits_{w\left(  P\right)  }^{w\left(
Q\right)  }A\right)  \delta_{P,Q}\prod_{i\in\left\{  1,2,\ldots,n\right\}
\setminus P}d_{i}\nonumber\\
&  =\underbrace{\left(  -1\right)  ^{\sum P+\sum P}}%
_{\substack{=1\\\text{(since }\sum P+\sum P=2\sum P\text{ is even)}}%
}\det\left(  \operatorname*{sub}\nolimits_{w\left(  P\right)  }^{w\left(
P\right)  }A\right)  \underbrace{\delta_{P,P}}_{=1}\prod_{i\in\left\{
1,2,\ldots,n\right\}  \setminus P}d_{i}\nonumber\\
&  \ \ \ \ \ \ \ \ \ \ +\sum_{\substack{Q\subseteq\left\{  1,2,\ldots
,n\right\}  ;\\\left\vert P\right\vert =\left\vert Q\right\vert ;\\Q\neq
P}}\left(  -1\right)  ^{\sum P+\sum Q}\det\left(  \operatorname*{sub}%
\nolimits_{w\left(  P\right)  }^{w\left(  Q\right)  }A\right)
\underbrace{\delta_{P,Q}}_{\substack{=0\\\text{(since }P\neq Q\\\text{(since
}Q\neq P\text{))}}}\prod_{i\in\left\{  1,2,\ldots,n\right\}  \setminus P}%
d_{i}\nonumber\\
&  \ \ \ \ \ \ \ \ \ \ \left(
\begin{array}
[c]{c}%
\text{here, we have split off the addend for }Q=P\text{ from the sum}\\
\text{(since }P\text{ is a subset }Q\text{ of }\left\{  1,2,\ldots,n\right\}
\text{ satisfying }\left\vert P\right\vert =\left\vert Q\right\vert \text{)}%
\end{array}
\right) \nonumber\\
&  =\det\left(  \operatorname*{sub}\nolimits_{w\left(  P\right)  }^{w\left(
P\right)  }A\right)  \prod_{i\in\left\{  1,2,\ldots,n\right\}  \setminus
P}d_{i}\nonumber\\
&  \ \ \ \ \ \ \ \ \ \ +\underbrace{\sum_{\substack{Q\subseteq\left\{
1,2,\ldots,n\right\}  ;\\\left\vert P\right\vert =\left\vert Q\right\vert
;\\Q\neq P}}\left(  -1\right)  ^{\sum P+\sum Q}\det\left(  \operatorname*{sub}%
\nolimits_{w\left(  P\right)  }^{w\left(  Q\right)  }A\right)  0\prod
_{i\in\left\{  1,2,\ldots,n\right\}  \setminus P}d_{i}}_{=0}\nonumber\\
&  =\det\left(  \operatorname*{sub}\nolimits_{w\left(  P\right)  }^{w\left(
P\right)  }A\right)  \prod_{i\in\left\{  1,2,\ldots,n\right\}  \setminus
P}d_{i}. \label{pf.cor.det(A+D).short.6}%
\end{align}

\end{vershort}

\begin{verlong}
Now, let $P$ be a subset of $\left\{  1,2,\ldots,n\right\}  $. Then,%
\begin{align}
&  \sum_{\substack{Q\subseteq\left\{  1,2,\ldots,n\right\}  ;\\\left\vert
P\right\vert =\left\vert Q\right\vert }}\left(  -1\right)  ^{\sum P+\sum
Q}\det\left(  \operatorname*{sub}\nolimits_{w\left(  P\right)  }^{w\left(
Q\right)  }A\right)  \underbrace{\det\left(  \operatorname*{sub}%
\nolimits_{w\left(  \widetilde{P}\right)  }^{w\left(  \widetilde{Q}\right)
}D\right)  }_{\substack{=\delta_{P,Q}\prod_{i\in\left\{  1,2,\ldots,n\right\}
\setminus P}d_{i}\\\text{(by (\ref{pf.cor.det(A+D).1}))}}}\nonumber\\
&  =\sum_{\substack{Q\subseteq\left\{  1,2,\ldots,n\right\}  ;\\\left\vert
P\right\vert =\left\vert Q\right\vert }}\left(  -1\right)  ^{\sum P+\sum
Q}\det\left(  \operatorname*{sub}\nolimits_{w\left(  P\right)  }^{w\left(
Q\right)  }A\right)  \delta_{P,Q}\prod_{i\in\left\{  1,2,\ldots,n\right\}
\setminus P}d_{i}\nonumber\\
&  =\sum_{\substack{Q\subseteq\left\{  1,2,\ldots,n\right\}  ;\\\left\vert
P\right\vert =\left\vert Q\right\vert ;\\Q=P}}\left(  -1\right)  ^{\sum P+\sum
Q}\det\left(  \operatorname*{sub}\nolimits_{w\left(  P\right)  }^{w\left(
Q\right)  }A\right)  \underbrace{\delta_{P,Q}}_{\substack{=1\\\text{(since
}P=Q\\\text{(since }Q=P\text{))}}}\prod_{i\in\left\{  1,2,\ldots,n\right\}
\setminus P}d_{i}\nonumber\\
&  \ \ \ \ \ \ \ \ \ \ +\sum_{\substack{Q\subseteq\left\{  1,2,\ldots
,n\right\}  ;\\\left\vert P\right\vert =\left\vert Q\right\vert ;\\Q\neq
P}}\left(  -1\right)  ^{\sum P+\sum Q}\det\left(  \operatorname*{sub}%
\nolimits_{w\left(  P\right)  }^{w\left(  Q\right)  }A\right)
\underbrace{\delta_{P,Q}}_{\substack{=0\\\text{(since }P\neq Q\\\text{(since
}Q\neq P\text{))}}}\prod_{i\in\left\{  1,2,\ldots,n\right\}  \setminus P}%
d_{i}\nonumber\\
&  \ \ \ \ \ \ \ \ \ \ \left(
\begin{array}
[c]{c}%
\text{since every subset }Q\text{ of }\left\{  1,2,\ldots,n\right\}  \text{
satisfying }\left\vert P\right\vert =\left\vert Q\right\vert \\
\text{satisfies either }Q=P\text{ or }Q\neq P\text{ (but not both)}%
\end{array}
\right) \nonumber\\
&  =\sum_{\substack{Q\subseteq\left\{  1,2,\ldots,n\right\}  ;\\\left\vert
P\right\vert =\left\vert Q\right\vert ;\\Q=P}}\left(  -1\right)  ^{\sum P+\sum
Q}\det\left(  \operatorname*{sub}\nolimits_{w\left(  P\right)  }^{w\left(
Q\right)  }A\right)  \prod_{i\in\left\{  1,2,\ldots,n\right\}  \setminus
P}d_{i}\nonumber\\
&  \ \ \ \ \ \ \ \ \ \ +\underbrace{\sum_{\substack{Q\subseteq\left\{
1,2,\ldots,n\right\}  ;\\\left\vert P\right\vert =\left\vert Q\right\vert
;\\Q\neq P}}\left(  -1\right)  ^{\sum P+\sum Q}\det\left(  \operatorname*{sub}%
\nolimits_{w\left(  P\right)  }^{w\left(  Q\right)  }A\right)  0\prod
_{i\in\left\{  1,2,\ldots,n\right\}  \setminus P}d_{i}}_{=0}\nonumber\\
&  =\sum_{\substack{Q\subseteq\left\{  1,2,\ldots,n\right\}  ;\\\left\vert
P\right\vert =\left\vert Q\right\vert ;\\Q=P}}\left(  -1\right)  ^{\sum P+\sum
Q}\det\left(  \operatorname*{sub}\nolimits_{w\left(  P\right)  }^{w\left(
Q\right)  }A\right)  \prod_{i\in\left\{  1,2,\ldots,n\right\}  \setminus
P}d_{i}. \label{pf.cor.det(A+D).4}%
\end{align}
But $P$ is a subset $Q$ of $\left\{  1,2,\ldots,n\right\}  $ satisfying
$\left\vert P\right\vert =\left\vert Q\right\vert $ and $Q=P$ (because
$\left\vert P\right\vert =\left\vert P\right\vert $ and $P=P$). Also, $P$ is
clearly the only such subset $Q$ (because any such subset $Q$ must satisfy
$Q=P$). Hence, there exists exactly one subset $Q$ of $\left\{  1,2,\ldots
,n\right\}  $ satisfying $\left\vert P\right\vert =\left\vert Q\right\vert $
and $Q=P$: namely, the subset $P$. Thus, the sum
\[
\sum_{\substack{Q\subseteq\left\{  1,2,\ldots,n\right\}  ;\\\left\vert
P\right\vert =\left\vert Q\right\vert ;\\Q=P}}\left(  -1\right)  ^{\sum P+\sum
Q}\det\left(  \operatorname*{sub}\nolimits_{w\left(  P\right)  }^{w\left(
Q\right)  }A\right)  \prod_{i\in\left\{  1,2,\ldots,n\right\}  \setminus
P}d_{i}%
\]
has only one addend: namely, the addend for $Q=P$. Thus, this sum rewrites as
follows:%
\begin{align*}
&  \sum_{\substack{Q\subseteq\left\{  1,2,\ldots,n\right\}  ;\\\left\vert
P\right\vert =\left\vert Q\right\vert ;\\Q=P}}\left(  -1\right)  ^{\sum P+\sum
Q}\det\left(  \operatorname*{sub}\nolimits_{w\left(  P\right)  }^{w\left(
Q\right)  }A\right)  \prod_{i\in\left\{  1,2,\ldots,n\right\}  \setminus
P}d_{i}\\
&  =\underbrace{\left(  -1\right)  ^{\sum P+\sum P}}%
_{\substack{=1\\\text{(since }\sum P+\sum P=2\sum P\text{ is even)}}%
}\det\left(  \operatorname*{sub}\nolimits_{w\left(  P\right)  }^{w\left(
P\right)  }A\right)  \prod_{i\in\left\{  1,2,\ldots,n\right\}  \setminus
P}d_{i}\\
&  =\det\left(  \operatorname*{sub}\nolimits_{w\left(  P\right)  }^{w\left(
P\right)  }A\right)  \prod_{i\in\left\{  1,2,\ldots,n\right\}  \setminus
P}d_{i}.
\end{align*}
Now, (\ref{pf.cor.det(A+D).4}) becomes%
\begin{align}
&  \sum_{\substack{Q\subseteq\left\{  1,2,\ldots,n\right\}  ;\\\left\vert
P\right\vert =\left\vert Q\right\vert }}\left(  -1\right)  ^{\sum P+\sum
Q}\det\left(  \operatorname*{sub}\nolimits_{w\left(  P\right)  }^{w\left(
Q\right)  }A\right)  \det\left(  \operatorname*{sub}\nolimits_{w\left(
\widetilde{P}\right)  }^{w\left(  \widetilde{Q}\right)  }D\right) \nonumber\\
&  =\sum_{\substack{Q\subseteq\left\{  1,2,\ldots,n\right\}  ;\\\left\vert
P\right\vert =\left\vert Q\right\vert ;\\Q=P}}\left(  -1\right)  ^{\sum P+\sum
Q}\det\left(  \operatorname*{sub}\nolimits_{w\left(  P\right)  }^{w\left(
Q\right)  }A\right)  \prod_{i\in\left\{  1,2,\ldots,n\right\}  \setminus
P}d_{i}\nonumber\\
&  =\det\left(  \operatorname*{sub}\nolimits_{w\left(  P\right)  }^{w\left(
P\right)  }A\right)  \prod_{i\in\left\{  1,2,\ldots,n\right\}  \setminus
P}d_{i}. \label{pf.cor.det(A+D).6}%
\end{align}

Now, forget that we fixed $P$. We thus have proven that
(\ref{pf.cor.det(A+D).6}) holds for every subset $P$ of $\left\{
1,2,\ldots,n\right\}  $.
\end{verlong}

\begin{vershort}
Now, Theorem \ref{thm.det(A+B)} (applied to $B=D$) yields%
\begin{align*}
&  \det\left(  A+D\right) \\
&  =\sum_{P\subseteq\left\{  1,2,\ldots,n\right\}  }\underbrace{\sum
_{\substack{Q\subseteq\left\{  1,2,\ldots,n\right\}  ;\\\left\vert
P\right\vert =\left\vert Q\right\vert }}\left(  -1\right)  ^{\sum P+\sum
Q}\det\left(  \operatorname*{sub}\nolimits_{w\left(  P\right)  }^{w\left(
Q\right)  }A\right)  \det\left(  \operatorname*{sub}\nolimits_{w\left(
\widetilde{P}\right)  }^{w\left(  \widetilde{Q}\right)  }D\right)
}_{\substack{=\det\left(  \operatorname*{sub}\nolimits_{w\left(  P\right)
}^{w\left(  P\right)  }A\right)  \prod_{i\in\left\{  1,2,\ldots,n\right\}
\setminus P}d_{i}\\\text{(by (\ref{pf.cor.det(A+D).short.6}))}}}\\
&  =\sum_{P\subseteq\left\{  1,2,\ldots,n\right\}  }\det\left(
\operatorname*{sub}\nolimits_{w\left(  P\right)  }^{w\left(  P\right)
}A\right)  \prod_{i\in\left\{  1,2,\ldots,n\right\}  \setminus P}d_{i}.
\end{align*}
This proves Corollary \ref{cor.det(A+D)}. \qedhere

\end{vershort}

\begin{verlong}
Now, Theorem \ref{thm.det(A+B)} (applied to $B=D$) yields%
\begin{align*}
&  \det\left(  A+D\right) \\
&  =\sum_{P\subseteq\left\{  1,2,\ldots,n\right\}  }\underbrace{\sum
_{\substack{Q\subseteq\left\{  1,2,\ldots,n\right\}  ;\\\left\vert
P\right\vert =\left\vert Q\right\vert }}\left(  -1\right)  ^{\sum P+\sum
Q}\det\left(  \operatorname*{sub}\nolimits_{w\left(  P\right)  }^{w\left(
Q\right)  }A\right)  \det\left(  \operatorname*{sub}\nolimits_{w\left(
\widetilde{P}\right)  }^{w\left(  \widetilde{Q}\right)  }D\right)
}_{\substack{=\det\left(  \operatorname*{sub}\nolimits_{w\left(  P\right)
}^{w\left(  P\right)  }A\right)  \prod_{i\in\left\{  1,2,\ldots,n\right\}
\setminus P}d_{i}\\\text{(by (\ref{pf.cor.det(A+D).6}))}}}\\
&  =\sum_{P\subseteq\left\{  1,2,\ldots,n\right\}  }\det\left(
\operatorname*{sub}\nolimits_{w\left(  P\right)  }^{w\left(  P\right)
}A\right)  \prod_{i\in\left\{  1,2,\ldots,n\right\}  \setminus P}d_{i}.
\end{align*}
This proves Corollary \ref{cor.det(A+D)}.
\end{verlong}
\end{proof}

\begin{proof}
[Proof of Corollary \ref{cor.det(A+X)}.]For every two objects $i$ and $j$,
define $\delta_{i,j}\in\mathbb{K}$ by $\delta_{i,j}=%
\begin{cases}
1, & \text{if }i=j;\\
0, & \text{if }i\neq j
\end{cases}
$. Then, $I_{n}=\left(  \delta_{i,j}\right)  _{1\leq i\leq n,\ 1\leq j\leq n}$
(by the definition of $I_{n}$). Now,%
\[
x\underbrace{I_{n}}_{=\left(  \delta_{i,j}\right)  _{1\leq i\leq n,\ 1\leq
j\leq n}}=x\left(  \delta_{i,j}\right)  _{1\leq i\leq n,\ 1\leq j\leq
n}=\left(  x\delta_{i,j}\right)  _{1\leq i\leq n,\ 1\leq j\leq n}%
\]
(by the definition of $x\left(  \delta_{i,j}\right)  _{1\leq i\leq n,\ 1\leq
j\leq n}$). Hence, Corollary \ref{cor.det(A+D)} (applied to $d_{i}=x$ and
$D=xI_{n}$) yields%
\begin{align*}
\det\left(  A+xI_{n}\right)   &  =\sum_{P\subseteq\left\{  1,2,\ldots
,n\right\}  }\det\left(  \operatorname*{sub}\nolimits_{w\left(  P\right)
}^{w\left(  P\right)  }A\right)  \underbrace{\prod_{i\in\left\{
1,2,\ldots,n\right\}  \setminus P}x}_{\substack{=x^{\left\vert \left\{
1,2,\ldots,n\right\}  \setminus P\right\vert }=x^{\left\vert \left\{
1,2,\ldots,n\right\}  \right\vert -\left\vert P\right\vert }\\\text{(since
}\left\vert \left\{  1,2,\ldots,n\right\}  \setminus P\right\vert =\left\vert
\left\{  1,2,\ldots,n\right\}  \right\vert -\left\vert P\right\vert
\\\text{(since }P\subseteq\left\{  1,2,\ldots,n\right\}  \text{))}}}\\
&  =\sum_{P\subseteq\left\{  1,2,\ldots,n\right\}  }\det\left(
\operatorname*{sub}\nolimits_{w\left(  P\right)  }^{w\left(  P\right)
}A\right)  \underbrace{x^{\left\vert \left\{  1,2,\ldots,n\right\}
\right\vert -\left\vert P\right\vert }}_{\substack{=x^{n-\left\vert
P\right\vert }\\\text{(since }\left\vert \left\{  1,2,\ldots,n\right\}
\right\vert =n\text{)}}}\\
&  =\sum_{P\subseteq\left\{  1,2,\ldots,n\right\}  }\det\left(
\operatorname*{sub}\nolimits_{w\left(  P\right)  }^{w\left(  P\right)
}A\right)  x^{n-\left\vert P\right\vert }.
\end{align*}
This proves (\ref{eq.cor.det(A+X).1}).

Furthermore, every subset $P$ of $\left\{  1,2,\ldots,n\right\}  $ satisfies
$\left\vert P\right\vert \in\left\{  0,1,\ldots,n\right\}  $%
\ \ \ \ \footnote{\textit{Proof.} Let $P$ be a subset of $\left\{
1,2,\ldots,n\right\}  $. Then, $P$ is a finite set (since $\left\{
1,2,\ldots,n\right\}  $ is a finite set), so that $\left\vert P\right\vert
\in\mathbb{N}$. Also, $P$ is a subset of $\left\{  1,2,\ldots,n\right\}  $,
and thus we have $\left\vert P\right\vert \leq\left\vert \left\{
1,2,\ldots,n\right\}  \right\vert =n$. Combined with $\left\vert P\right\vert
\in\mathbb{N}$, this yields $\left\vert P\right\vert \in\left\{
0,1,\ldots,n\right\}  $. Qed.}. Now,
\begin{align*}
\det\left(  A+xI_{n}\right)   &  =\underbrace{\sum_{P\subseteq\left\{
1,2,\ldots,n\right\}  }}_{\substack{=\sum_{k\in\left\{  0,1,\ldots,n\right\}
}\sum_{\substack{P\subseteq\left\{  1,2,\ldots,n\right\}  ;\\\left\vert
P\right\vert =k}}\\\text{(since every subset }P\text{ of }\left\{
1,2,\ldots,n\right\}  \\\text{satisfies }\left\vert P\right\vert \in\left\{
0,1,\ldots,n\right\}  \text{)}}}\det\left(  \operatorname*{sub}%
\nolimits_{w\left(  P\right)  }^{w\left(  P\right)  }A\right)  x^{n-\left\vert
P\right\vert }\\
&  =\underbrace{\sum_{k\in\left\{  0,1,\ldots,n\right\}  }}_{=\sum_{k=0}^{n}%
}\sum_{\substack{P\subseteq\left\{  1,2,\ldots,n\right\}  ;\\\left\vert
P\right\vert =k}}\det\left(  \operatorname*{sub}\nolimits_{w\left(  P\right)
}^{w\left(  P\right)  }A\right)  \underbrace{x^{n-\left\vert P\right\vert }%
}_{\substack{=x^{n-k}\\\text{(since }\left\vert P\right\vert =k\text{)}}}\\
&  =\sum_{k=0}^{n}\sum_{\substack{P\subseteq\left\{  1,2,\ldots,n\right\}
;\\\left\vert P\right\vert =k}}\det\left(  \operatorname*{sub}%
\nolimits_{w\left(  P\right)  }^{w\left(  P\right)  }A\right)  x^{n-k}\\
&  =\sum_{k=0}^{n}\sum_{\substack{P\subseteq\left\{  1,2,\ldots,n\right\}
;\\\left\vert P\right\vert =n-k}}\det\left(  \operatorname*{sub}%
\nolimits_{w\left(  P\right)  }^{w\left(  P\right)  }A\right)
\underbrace{x^{n-\left(  n-k\right)  }}_{\substack{=x^{k}\\\text{(since
}n-\left(  n-k\right)  =k\text{)}}}\\
&  \ \ \ \ \ \ \ \ \ \ \left(  \text{here, we have substituted }n-k\text{ for
}k\text{ in the outer sum}\right) \\
&  =\sum_{k=0}^{n}\sum_{\substack{P\subseteq\left\{  1,2,\ldots,n\right\}
;\\\left\vert P\right\vert =n-k}}\det\left(  \operatorname*{sub}%
\nolimits_{w\left(  P\right)  }^{w\left(  P\right)  }A\right)  x^{k}\\
&  =\sum_{k=0}^{n}\left(  \sum_{\substack{P\subseteq\left\{  1,2,\ldots
,n\right\}  ;\\\left\vert P\right\vert =n-k}}\det\left(  \operatorname*{sub}%
\nolimits_{w\left(  P\right)  }^{w\left(  P\right)  }A\right)  \right)  x^{k}.
\end{align*}
This proves (\ref{eq.cor.det(A+X).2}). The proof of Corollary
\ref{cor.det(A+X)} is now complete.
\end{proof}

\begin{proof}
[Solution to Exercise \ref{exe.det(A+B).diag}.]We have proven Lemma
\ref{lem.diag.minors}, Corollary \ref{cor.det(A+D)} and Corollary
\ref{cor.det(A+X)}. Exercise \ref{exe.det(A+B).diag} is thus solved.
\end{proof}

\subsection{\label{sect.sol.noncomm.polarization}Solution to Exercise
\ref{exe.noncomm.polarization}}

Let us first prepare some auxiliary results and notations. Throughout Section
\ref{sect.sol.noncomm.polarization}, we shall be using the following conventions:

\begin{itemize}
\item For every $n\in\mathbb{N}$, let $\left[  n\right]  $ denote the set
$\left\{  1,2,\ldots,n\right\}  $.

\item For every $n\in\mathbb{N}$, the sign $\sum_{I\subseteq\left[  n\right]
}$ shall mean $\sum_{I\in\mathcal{P}\left(  \left[  n\right]  \right)  }$,
where $\mathcal{P}\left(  \left[  n\right]  \right)  $ denotes the powerset of
$\left[  n\right]  $.

\item We shall use the Iverson bracket notation introduced in Definition
\ref{def.iverson}.

\item If $a_{p},a_{p+1},\ldots,a_{q}$ are some elements of a noncommutative
ring $\mathbb{L}$ (with $p$ and $q$ being integers satisfying $p\leq q+1$),
then $\prod_{i=p}^{q}a_{i}$ is defined to be the product $a_{p}a_{p+1}\cdots
a_{q}\in\mathbb{L}$\ \ \ \ \footnote{Of course, in order for this definition
to be fully rigorous, we need to specify how the product $a_{p}a_{p+1}\cdots
a_{q}$ is defined. There are two possible definitions for this product; both
of them proceed by recursion on $q-p$ (with the recursion base being the case
when $q-p=-1$).
\par
\begin{itemize}
\item The first definition sets $a_{p}a_{p+1}\cdots a_{q}$ to be $1$ (the
unity of $\mathbb{L}$) when $q-p=-1$ (this is the case of an empty product),
and sets $a_{p}a_{p+1}\cdots a_{q}$ to be $\left(  a_{p}a_{p+1}\cdots
a_{q-1}\right)  a_{q}$ when $q-p\geq0$.
\par
\item The second definition sets $a_{p}a_{p+1}\cdots a_{q}$ to be $1$ (the
unity of $\mathbb{L}$) when $q-p=-1$ (this is the case of an empty product),
and sets $a_{p}a_{p+1}\cdots a_{q}$ to be $a_{p}\left(  a_{p+1}a_{p+2}\cdots
a_{q}\right)  $ when $q-p\geq0$.
\end{itemize}
\par
Fortunately, the two definitions define exactly the same notion of product;
this is not hard to prove (it follows from a fact called \textquotedblleft
general associativity\textquotedblright, which is analogous to Proposition
\ref{prop.ind.gen-ass-maps.Ceq-cp} but involves multiplication of elements of
$\mathbb{L}$ instead of composition of maps), but it is not entirely
obvious.}. This definition of $\prod_{i=p}^{q}a_{i}$ extends the classical
definition of $\prod_{i=p}^{q}a_{i}$ in the case when $a_{p},a_{p+1}%
,\ldots,a_{q}$ are elements of a commutative ring.

(Notice that we have thus defined products of the form $\prod_{i=p}^{q}a_{i}$
only. In contrast, products of the form $\prod_{i\in I}a_{i}$ for an arbitrary
finite set $I$ are \textbf{not} defined in a noncommutative ring, because it
is not clear in what order their factors are to be multiplied.)
\end{itemize}

A few elementary remarks about noncommutative rings will be useful:

\begin{itemize}
\item If $p$ and $q$ are integers satisfying $p\leq q$, and if $a_{p}%
,a_{p+1},\ldots,a_{q}$ are some elements of a noncommutative ring $\mathbb{L}%
$, then%
\[
\prod_{i=p}^{q}a_{i}=\left(  \prod_{i=p}^{q-1}a_{i}\right)  a_{q}.
\]
In other words, we can always split off the factor for $i=q$ from a product of
the form $\prod_{i=p}^{q}a_{i}$ as long as we have $p\leq q$.

\item When $\mathbb{L}$ is a noncommutative ring, the following two properties
of the $\sum$ sign hold\footnote{These are analogues of the \textquotedblleft
Factoring out\textquotedblright\ property that we have seen in Section
\ref{sect.sums-repetitorium}.}:

\begin{itemize}
\item \underline{\textbf{Factoring out on the left:}} Let $S$ be a finite set.
For every $s\in S$, let $a_{s}$ be an element of $\mathbb{L}$. Also, let
$\lambda$ be an element of $\mathbb{L}$. Then,%
\begin{equation}
\sum_{s\in S}\lambda a_{s}=\lambda\sum_{s\in S}a_{s}.
\label{eq.noncomm.sum.linear2}%
\end{equation}

\item \underline{\textbf{Factoring out on the right:}} Let $S$ be a finite
set. For every $s\in S$, let $a_{s}$ be an element of $\mathbb{L}$. Also, let
$\lambda$ be an element of $\mathbb{L}$. Then,%
\begin{equation}
\sum_{s\in S}a_{s}\lambda=\left(  \sum_{s\in S}a_{s}\right)  \lambda.
\label{eq.noncomm.sum.linear2r}%
\end{equation}

\end{itemize}
\end{itemize}

Next, we state a generalization of Lemma \ref{lem.prodrule2} to the case of a
noncommutative ring $\mathbb{K}$:

\begin{lemma}
\label{lem.noncomm.prodrule2}Let $\mathbb{K}$ be a noncommutative ring.

Let $n\in\mathbb{N}$ and $m\in\mathbb{N}$. For every $i\in\left[  n\right]  $,
let $p_{i,1},p_{i,2},\ldots,p_{i,m}$ be $m$ elements of $\mathbb{K}$. Then,%
\[
\prod_{i=1}^{n}\sum_{k=1}^{m}p_{i,k}=\sum_{\kappa:\left[  n\right]
\rightarrow\left[  m\right]  }\prod_{i=1}^{n}p_{i,\kappa\left(  i\right)  }.
\]

\end{lemma}

\begin{proof}
[Proof of Lemma \ref{lem.noncomm.prodrule2}.]The proof of Lemma
\ref{lem.prodrule2} given above (which includes the proof of Lemma
\ref{lem.prodrule} in the solution to Exercise \ref{exe.prodrule}) can be
reused as a proof of Lemma \ref{lem.noncomm.prodrule2}, provided that we
replace our proof of Lemma \ref{lem.prodrule.S.n=2} (which made use of the
commutativity of $\mathbb{K}$) by the following alternative proof (which does
not require $\mathbb{K}$ to be commutative):

\begin{proof}
[Second proof of Lemma \ref{lem.prodrule.S.n=2}.]The equality
(\ref{eq.sum.fubini}) remains valid if $\mathbb{A}$ is replaced by
$\mathbb{K}$ throughout it. It is in this variant that it will be used in the
following proof.

For every $\lambda\in\mathbb{K}$, we have%
\begin{align}
\sum_{y\in Y}\lambda r_{y}  &  =\sum_{s\in Y}\lambda r_{s}%
\ \ \ \ \ \ \ \ \ \ \left(  \text{here, we renamed the summation index
}y\text{ as }s\right) \nonumber\\
&  =\lambda\sum_{s\in Y}r_{s}\ \ \ \ \ \ \ \ \ \ \left(
\begin{array}
[c]{c}%
\text{by (\ref{eq.noncomm.sum.linear2}) (applied to }Y\text{, }r_{s}\text{ and
}\mathbb{K}\\
\text{instead of }S\text{, }a_{s}\text{ and }\mathbb{L}\text{)}%
\end{array}
\right) \nonumber\\
&  =\lambda\sum_{y\in Y}r_{y}\ \ \ \ \ \ \ \ \ \ \left(  \text{here, we
renamed the summation index }s\text{ as }y\right)  .
\label{pf.lem.noncomm.prodrule.S.n=2.1}%
\end{align}
For every $\lambda\in\mathbb{K}$, we have%
\begin{align}
\sum_{x\in X}q_{x}\lambda &  =\sum_{s\in X}q_{s}\lambda
\ \ \ \ \ \ \ \ \ \ \left(  \text{here, we renamed the summation index
}x\text{ as }s\right) \nonumber\\
&  =\left(  \sum_{s\in X}q_{s}\right)  \lambda\ \ \ \ \ \ \ \ \ \ \left(
\begin{array}
[c]{c}%
\text{by (\ref{eq.noncomm.sum.linear2r}) (applied to }X\text{, }q_{s}\text{
and }\mathbb{K}\\
\text{instead of }S\text{, }a_{s}\text{ and }\mathbb{L}\text{)}%
\end{array}
\right) \nonumber\\
&  =\left(  \sum_{x\in X}q_{x}\right)  \lambda\ \ \ \ \ \ \ \ \ \ \left(
\text{here, we renamed the summation index }s\text{ as }x\right)  .
\label{pf.lem.noncomm.prodrule.S.n=2.2}%
\end{align}

From (\ref{eq.sum.fubini}) (applied to $\mathbb{K}$ and $q_{x}r_{y}$ instead
of $\mathbb{A}$ and $a_{\left(  x,y\right)  }$), we obtain%
\[
\sum_{x\in X}\sum_{y\in Y}q_{x}r_{y}=\sum_{\left(  x,y\right)  \in X\times
Y}q_{x}r_{y}=\sum_{y\in Y}\sum_{x\in X}q_{x}r_{y}.
\]
Hence,%
\begin{align*}
\sum_{\left(  x,y\right)  \in X\times Y}q_{x}r_{y}  &  =\sum_{x\in
X}\underbrace{\sum_{y\in Y}q_{x}r_{y}}_{\substack{=q_{x}\sum_{y\in Y}%
r_{y}\\\text{(by (\ref{pf.lem.noncomm.prodrule.S.n=2.1}) (applied to }%
\lambda=q_{x}\text{))}}}=\sum_{x\in X}q_{x}\sum_{y\in Y}r_{y}\\
&  =\left(  \sum_{x\in X}q_{x}\right)  \left(  \sum_{y\in Y}r_{y}\right)
\end{align*}
(by (\ref{pf.lem.noncomm.prodrule.S.n=2.2}) (applied to $\lambda=\sum_{y\in
Y}r_{y}$)). This proves Lemma \ref{lem.prodrule.S.n=2}.
\end{proof}

Thus, Lemma \ref{lem.noncomm.prodrule2} is proven.
\end{proof}

Now, let us prove some further lemmas:

\begin{lemma}
\label{lem.sol.noncomm.polarization.1.7}Let $G$ be a finite set. Let $I$ be a
subset of $G$. Let $m\in\mathbb{N}$. Let $f:\left[  m\right]  \rightarrow G$
be any map. Then,%
\[
\prod_{i=1}^{m}\left[  f\left(  i\right)  \in I\right]  =\left[  f\left(
\left[  m\right]  \right)  \subseteq I\right]  .
\]

\end{lemma}

\begin{vershort}
\begin{proof}
[Proof of Lemma \ref{lem.sol.noncomm.polarization.1.7}.]We are in one of the
following two cases:

\textit{Case 1:} We have $f\left(  \left[  m\right]  \right)  \subseteq I$.

\textit{Case 2:} We don't have $f\left(  \left[  m\right]  \right)  \subseteq
I$.

Let us first consider Case 1. In this case, we have $f\left(  \left[
m\right]  \right)  \subseteq I$. Thus, $\left[  f\left(  \left[  m\right]
\right)  \subseteq I\right]  =1$. But each $i\in\left\{  1,2,\ldots,m\right\}
$ satisfies $\left[  f\left(  i\right)  \in I\right]  =1$ (since $f\left(
\underbrace{i}_{\in\left[  m\right]  }\right)  \in f\left(  \left[  m\right]
\right)  \subseteq I$). Hence, $\prod_{i=1}^{m}\underbrace{\left[  f\left(
i\right)  \in I\right]  }_{=1}=\prod_{i=1}^{m}1=1$. Comparing this with
$\left[  f\left(  \left[  m\right]  \right)  \subseteq I\right]  =1$, we
obtain $\prod_{i=1}^{m}\left[  f\left(  i\right)  \in I\right]  =\left[
f\left(  \left[  m\right]  \right)  \subseteq I\right]  $. Hence, Lemma
\ref{lem.sol.noncomm.polarization.1.7} is proven in Case 1.

Let us now consider Case 2. In this case, we don't have $f\left(  \left[
m\right]  \right)  \subseteq I$. In other words, not every $p\in\left[
m\right]  $ satisfies $f\left(  p\right)  \in I$. Thus, there exists some
$p\in\left[  m\right]  $ such that $f\left(  p\right)  \notin I$. Consider
this $p$. We have $\left[  f\left(  p\right)  \in I\right]  =0$ (since
$f\left(  p\right)  \notin I$). Thus, at least one factor of the product
$\prod_{i=1}^{m}\left[  f\left(  i\right)  \in I\right]  $ equals $0$ (namely,
the factor for $i=p$ is $\left[  f\left(  p\right)  \in I\right]  =0$).
Consequently, the whole product $\prod_{i=1}^{m}\left[  f\left(  i\right)  \in
I\right]  $ equals $0$ (because if at least one factor of a product equals
$0$, then the whole product equals $0$). In other words, $\prod_{i=1}%
^{m}\left[  f\left(  i\right)  \in I\right]  =0$. Comparing this with $\left[
f\left(  \left[  m\right]  \right)  \subseteq I\right]  =0$ (which holds
because we don't have $f\left(  \left[  m\right]  \right)  \subseteq I$), we
obtain $\prod_{i=1}^{m}\left[  f\left(  i\right)  \in I\right]  =\left[
f\left(  \left[  m\right]  \right)  \subseteq I\right]  $. Hence, Lemma
\ref{lem.sol.noncomm.polarization.1.7} is proven in Case 2.

We have now proven Lemma \ref{lem.sol.noncomm.polarization.1.7} in both Cases
1 and 2. This completes the proof of Lemma
\ref{lem.sol.noncomm.polarization.1.7}.
\end{proof}
\end{vershort}

\begin{verlong}
\begin{proof}
[Proof of Lemma \ref{lem.sol.noncomm.polarization.1.7}.]We are in one of the
following two cases:

\textit{Case 1:} We have $f\left(  \left[  m\right]  \right)  \subseteq I$.

\textit{Case 2:} We don't have $f\left(  \left[  m\right]  \right)  \subseteq
I$.

Let us first consider Case 1. In this case, we have $f\left(  \left[
m\right]  \right)  \subseteq I$. Thus, $\left[  f\left(  \left[  m\right]
\right)  \subseteq I\right]  =1$. But each $i\in\left\{  1,2,\ldots,m\right\}
$ satisfies $\left[  f\left(  i\right)  \in I\right]  =1$ (since $f\left(
\underbrace{i}_{\in\left\{  1,2,\ldots,m\right\}  =\left[  m\right]  }\right)
\in f\left(  \left[  m\right]  \right)  \subseteq I$). Hence, $\prod_{i=1}%
^{m}\underbrace{\left[  f\left(  i\right)  \in I\right]  }_{=1}=\prod
_{i=1}^{m}1=1$. Comparing this with $\left[  f\left(  \left[  m\right]
\right)  \subseteq I\right]  =1$, we obtain $\prod_{i=1}^{m}\left[  f\left(
i\right)  \in I\right]  =\left[  f\left(  \left[  m\right]  \right)  \subseteq
I\right]  $. Hence, Lemma \ref{lem.sol.noncomm.polarization.1.7} is proven in
Case 1.

Let us now consider Case 2. In this case, we don't have $f\left(  \left[
m\right]  \right)  \subseteq I$. Thus, $\left[  f\left(  \left[  m\right]
\right)  \subseteq I\right]  =0$. But we have $f\left(  \left[  m\right]
\right)  \not \subseteq I$ (since we don't have $f\left(  \left[  m\right]
\right)  \subseteq I$). Hence, there exists some $q\in f\left(  \left[
m\right]  \right)  $ such that $q\notin I$. Consider this $q$. We have $q\in
f\left(  \left[  m\right]  \right)  $. Thus, there exists some $p\in\left[
m\right]  $ such that $q=f\left(  p\right)  $. Consider this $p$. We have
$f\left(  p\right)  =q\notin I$. In other words, we don't have $f\left(
p\right)  \in I$. Hence, $\left[  f\left(  p\right)  \in I\right]  =0$. But we
have $p\in\left[  m\right]  =\left\{  1,2,\ldots,m\right\}  $. Hence, $\left[
f\left(  p\right)  \in I\right]  $ is a factor of the product $\prod_{i=1}%
^{m}\left[  f\left(  i\right)  \in I\right]  $ (namely, the factor for $i=p$).
Moreover, this factor equals $0$ (since $\left[  f\left(  p\right)  \in
I\right]  =0$). Hence, at least one factor of the product $\prod_{i=1}%
^{m}\left[  f\left(  i\right)  \in I\right]  $ equals $0$ (namely, the factor
$\left[  f\left(  p\right)  \in I\right]  $). Consequently, the whole product
$\prod_{i=1}^{m}\left[  f\left(  i\right)  \in I\right]  $ equals $0$ (because
if at least one factor of a product equals $0$, then the whole product equals
$0$). In other words, $\prod_{i=1}^{m}\left[  f\left(  i\right)  \in I\right]
=0$. Comparing this with $\left[  f\left(  \left[  m\right]  \right)
\subseteq I\right]  =0$, we obtain $\prod_{i=1}^{m}\left[  f\left(  i\right)
\in I\right]  =\left[  f\left(  \left[  m\right]  \right)  \subseteq I\right]
$. Hence, Lemma \ref{lem.sol.noncomm.polarization.1.7} is proven in Case 2.

We have now proven Lemma \ref{lem.sol.noncomm.polarization.1.7} in both Cases
1 and 2. Since these are the only possible cases, we thus have always proven
Lemma \ref{lem.sol.noncomm.polarization.1.7}.
\end{proof}
\end{verlong}

\begin{lemma}
\label{lem.sol.noncomm.polarization.1}Let $\mathbb{L}$ be a noncommutative
ring. Let $n\in\mathbb{N}$. Let $v_{1},v_{2},\ldots,v_{n}$ be $n$ elements of
$\mathbb{L}$. Let $I$ be a subset of $\left[  n\right]  $. Then,%
\[
\left(  \sum_{i\in I}v_{i}\right)  ^{m}=\sum_{\substack{f:\left[  m\right]
\rightarrow\left[  n\right]  ;\\f\left(  \left[  m\right]  \right)  \subseteq
I}}v_{f\left(  1\right)  }v_{f\left(  2\right)  }\cdots v_{f\left(  m\right)
}.
\]

\end{lemma}

\begin{proof}
[Proof of Lemma \ref{lem.sol.noncomm.polarization.1}.]First, observe that%
\begin{align}
&  \sum_{k=1}^{n}\left[  k\in I\right]  v_{k}\nonumber\\
&  =\underbrace{\sum_{i=1}^{n}}_{\substack{=\sum_{i\in\left\{  1,2,\ldots
,n\right\}  }=\sum_{i\in\left[  n\right]  }\\\text{(since }\left\{
1,2,\ldots,n\right\}  =\left[  n\right]  \text{)}}}\left[  i\in I\right]
v_{i}\ \ \ \ \ \ \ \ \ \ \left(
\begin{array}
[c]{c}%
\text{here, we have renamed the}\\
\text{summation index }k\text{ as }i
\end{array}
\right) \nonumber\\
&  =\sum_{i\in\left[  n\right]  }\left[  i\in I\right]  v_{i}\nonumber\\
&  =\underbrace{\sum_{\substack{i\in\left[  n\right]  ;\\i\in I}%
}}_{\substack{=\sum_{i\in I}\\\text{(since }I\subseteq\left[  n\right]
\text{)}}}\underbrace{\left[  i\in I\right]  }_{\substack{=1\\\text{(since
}i\in I\text{ is true)}}}v_{i}+\sum_{\substack{i\in\left[  n\right]
;\\i\notin I}}\underbrace{\left[  i\in I\right]  }_{\substack{=0\\\text{(since
}i\in I\text{ is false}\\\text{(because }i\notin I\text{))}}}v_{i}\nonumber\\
&  \ \ \ \ \ \ \ \ \ \ \left(  \text{since each }i\in\left[  n\right]  \text{
satisfies either }i\in I\text{ or }i\notin I\text{ (but not both)}\right)
\nonumber\\
&  =\sum_{i\in I}1v_{i}+\underbrace{\sum_{\substack{i\in\left[  n\right]
;\\i\notin I}}0v_{i}}_{=0}=\sum_{i\in I}1v_{i}=\sum_{i\in I}v_{i}.
\label{pf.lem.sol.noncomm.polarization.1.2}%
\end{align}

But $p^{m}=\prod_{i=1}^{m}p$ for each $p\in\mathbb{L}$. Applying this to
$p=\sum_{k=1}^{n}\left[  k\in I\right]  v_{k}$, we find%
\[
\left(  \sum_{k=1}^{n}\left[  k\in I\right]  v_{k}\right)  ^{m}=\prod
_{i=1}^{m}\sum_{k=1}^{n}\left[  k\in I\right]  v_{k}=\sum_{\kappa:\left[
m\right]  \rightarrow\left[  n\right]  }\prod_{i=1}^{m}\left(  \left[
\kappa\left(  i\right)  \in I\right]  v_{\kappa\left(  i\right)  }\right)
\]
(by Lemma \ref{lem.noncomm.prodrule2} (applied to $\mathbb{L}$, $m$, $n$ and
$\left[  k\in I\right]  v_{k}$ instead of $\mathbb{K}$, $n$, $m$ and $p_{i,k}%
$)). Now,%
\begin{align}
\left(  \underbrace{\sum_{i\in I}v_{i}}_{\substack{=\sum_{k=1}^{n}\left[  k\in
I\right]  v_{k}\\\text{(by (\ref{pf.lem.sol.noncomm.polarization.1.2}))}%
}}\right)  ^{m}  &  =\left(  \sum_{k=1}^{n}\left[  k\in I\right]
v_{k}\right)  ^{m}\nonumber\\
&  =\sum_{\kappa:\left[  m\right]  \rightarrow\left[  n\right]  }\prod
_{i=1}^{m}\left(  \left[  \kappa\left(  i\right)  \in I\right]  v_{\kappa
\left(  i\right)  }\right) \nonumber\\
&  =\sum_{f:\left[  m\right]  \rightarrow\left[  n\right]  }\prod_{i=1}%
^{m}\left(  \left[  f\left(  i\right)  \in I\right]  v_{f\left(  i\right)
}\right)  \label{pf.lem.sol.noncomm.polarization.1.4}%
\end{align}
(here, we have renamed the summation index $\kappa$ as $f$).

Now, fix any map $f:\left[  m\right]  \rightarrow\left[  n\right]  $. Then,%
\begin{align}
&  \prod_{i=1}^{m}\left(  \left[  f\left(  i\right)  \in I\right]  v_{f\left(
i\right)  }\right) \nonumber\\
&  =\left(  \left[  f\left(  1\right)  \in I\right]  v_{f\left(  1\right)
}\right)  \left(  \left[  f\left(  2\right)  \in I\right]  v_{f\left(
2\right)  }\right)  \cdots\left(  \left[  f\left(  m\right)  \in I\right]
v_{f\left(  m\right)  }\right)  . \label{pf.lem.sol.noncomm.polarization.1.6}%
\end{align}
The factors $\left[  f\left(  1\right)  \in I\right]  ,\left[  f\left(
2\right)  \in I\right]  ,\ldots,\left[  f\left(  m\right)  \in I\right]  $ on
the right hand side of this equality are integers, and therefore can be freely
moved within the product (even though $\mathbb{L}$ is not necessarily
commutative). In particular, we can move them to the front of the product.
Thus, we find%
\begin{align*}
&  \left(  \left[  f\left(  1\right)  \in I\right]  v_{f\left(  1\right)
}\right)  \left(  \left[  f\left(  2\right)  \in I\right]  v_{f\left(
2\right)  }\right)  \cdots\left(  \left[  f\left(  m\right)  \in I\right]
v_{f\left(  m\right)  }\right) \\
&  =\underbrace{\left(  \left[  f\left(  1\right)  \in I\right]  \left[
f\left(  2\right)  \in I\right]  \cdots\left[  f\left(  m\right)  \in
I\right]  \right)  }_{=\prod_{i=1}^{m}\left[  f\left(  i\right)  \in I\right]
}\left(  v_{f\left(  1\right)  }v_{f\left(  2\right)  }\cdots v_{f\left(
m\right)  }\right) \\
&  =\left(  \prod_{i=1}^{m}\left[  f\left(  i\right)  \in I\right]  \right)
\left(  v_{f\left(  1\right)  }v_{f\left(  2\right)  }\cdots v_{f\left(
m\right)  }\right)  .
\end{align*}
Hence, (\ref{pf.lem.sol.noncomm.polarization.1.6}) rewrites as%
\begin{equation}
\prod_{i=1}^{m}\left(  \left[  f\left(  i\right)  \in I\right]  v_{f\left(
i\right)  }\right)  =\left(  \prod_{i=1}^{m}\left[  f\left(  i\right)  \in
I\right]  \right)  \left(  v_{f\left(  1\right)  }v_{f\left(  2\right)
}\cdots v_{f\left(  m\right)  }\right)  .
\label{pf.lem.sol.noncomm.polarization.1.6b}%
\end{equation}

But Lemma \ref{lem.sol.noncomm.polarization.1.7} (applied to $G=\left[
n\right]  $) yields%
\begin{equation}
\prod_{i=1}^{m}\left[  f\left(  i\right)  \in I\right]  =\left[  f\left(
\left[  m\right]  \right)  \subseteq I\right]  .
\label{pf.lem.sol.noncomm.polarization.1.7n}%
\end{equation}

Hence, (\ref{pf.lem.sol.noncomm.polarization.1.6b}) becomes%
\begin{align}
\prod_{i=1}^{m}\left(  \left[  f\left(  i\right)  \in I\right]  v_{f\left(
i\right)  }\right)   &  =\underbrace{\left(  \prod_{i=1}^{m}\left[  f\left(
i\right)  \in I\right]  \right)  }_{\substack{=\left[  f\left(  \left[
m\right]  \right)  \subseteq I\right]  \\\text{(by
(\ref{pf.lem.sol.noncomm.polarization.1.7n}))}}}\left(  v_{f\left(  1\right)
}v_{f\left(  2\right)  }\cdots v_{f\left(  m\right)  }\right) \nonumber\\
&  =\left[  f\left(  \left[  m\right]  \right)  \subseteq I\right]  \left(
v_{f\left(  1\right)  }v_{f\left(  2\right)  }\cdots v_{f\left(  m\right)
}\right)  . \label{pf.lem.sol.noncomm.polarization.1.8}%
\end{align}

Now, forget that we fixed $f$. We thus have proven
(\ref{pf.lem.sol.noncomm.polarization.1.8}) for each map $f:\left[  m\right]
\rightarrow\left[  n\right]  $.

\begin{vershort}
Now, (\ref{pf.lem.sol.noncomm.polarization.1.4}) becomes%
\begin{align*}
\left(  \sum_{i\in I}v_{i}\right)  ^{m}  &  =\sum_{f:\left[  m\right]
\rightarrow\left[  n\right]  }\underbrace{\prod_{i=1}^{m}\left(  \left[
f\left(  i\right)  \in I\right]  v_{f\left(  i\right)  }\right)
}_{\substack{=\left[  f\left(  \left[  m\right]  \right)  \subseteq I\right]
\left(  v_{f\left(  1\right)  }v_{f\left(  2\right)  }\cdots v_{f\left(
m\right)  }\right)  \\\text{(by (\ref{pf.lem.sol.noncomm.polarization.1.8}))}%
}}\\
&  =\sum_{f:\left[  m\right]  \rightarrow\left[  n\right]  }\left[  f\left(
\left[  m\right]  \right)  \subseteq I\right]  \left(  v_{f\left(  1\right)
}v_{f\left(  2\right)  }\cdots v_{f\left(  m\right)  }\right) \\
&  =\sum_{\substack{f:\left[  m\right]  \rightarrow\left[  n\right]
;\\f\left(  \left[  m\right]  \right)  \subseteq I}}\underbrace{\left[
f\left(  \left[  m\right]  \right)  \subseteq I\right]  }%
_{\substack{=1\\\text{(since }f\left(  \left[  m\right]  \right)  \subseteq
I\text{)}}}\left(  v_{f\left(  1\right)  }v_{f\left(  2\right)  }\cdots
v_{f\left(  m\right)  }\right) \\
&  \ \ \ \ \ \ \ \ \ \ +\sum_{\substack{f:\left[  m\right]  \rightarrow\left[
n\right]  ;\\\text{not }f\left(  \left[  m\right]  \right)  \subseteq
I}}\ \ \underbrace{\left[  f\left(  \left[  m\right]  \right)  \subseteq
I\right]  }_{\substack{=0\\\text{(since we don't have }f\left(  \left[
m\right]  \right)  \subseteq I\text{)}}}\left(  v_{f\left(  1\right)
}v_{f\left(  2\right)  }\cdots v_{f\left(  m\right)  }\right) \\
&  =\sum_{\substack{f:\left[  m\right]  \rightarrow\left[  n\right]
;\\f\left(  \left[  m\right]  \right)  \subseteq I}}v_{f\left(  1\right)
}v_{f\left(  2\right)  }\cdots v_{f\left(  m\right)  }+\underbrace{\sum
_{\substack{f:\left[  m\right]  \rightarrow\left[  n\right]  ;\\\text{not
}f\left(  \left[  m\right]  \right)  \subseteq I}}0\left(  v_{f\left(
1\right)  }v_{f\left(  2\right)  }\cdots v_{f\left(  m\right)  }\right)
}_{=0}\\
&  =\sum_{\substack{f:\left[  m\right]  \rightarrow\left[  n\right]
;\\f\left(  \left[  m\right]  \right)  \subseteq I}}v_{f\left(  1\right)
}v_{f\left(  2\right)  }\cdots v_{f\left(  m\right)  }.
\end{align*}
This proves Lemma \ref{lem.sol.noncomm.polarization.1}. \qedhere

\end{vershort}

\begin{verlong}
Now, (\ref{pf.lem.sol.noncomm.polarization.1.4}) becomes%
\begin{align*}
\left(  \sum_{i\in I}v_{i}\right)  ^{m}  &  =\sum_{f:\left[  m\right]
\rightarrow\left[  n\right]  }\underbrace{\prod_{i=1}^{m}\left(  \left[
f\left(  i\right)  \in I\right]  v_{f\left(  i\right)  }\right)
}_{\substack{=\left[  f\left(  \left[  m\right]  \right)  \subseteq I\right]
\left(  v_{f\left(  1\right)  }v_{f\left(  2\right)  }\cdots v_{f\left(
m\right)  }\right)  \\\text{(by (\ref{pf.lem.sol.noncomm.polarization.1.8}))}%
}}\\
&  =\sum_{f:\left[  m\right]  \rightarrow\left[  n\right]  }\left[  f\left(
\left[  m\right]  \right)  \subseteq I\right]  \left(  v_{f\left(  1\right)
}v_{f\left(  2\right)  }\cdots v_{f\left(  m\right)  }\right) \\
&  =\sum_{\substack{f:\left[  m\right]  \rightarrow\left[  n\right]
;\\f\left(  \left[  m\right]  \right)  \subseteq I}}\underbrace{\left[
f\left(  \left[  m\right]  \right)  \subseteq I\right]  }%
_{\substack{=1\\\text{(since }f\left(  \left[  m\right]  \right)  \subseteq
I\text{)}}}\left(  v_{f\left(  1\right)  }v_{f\left(  2\right)  }\cdots
v_{f\left(  m\right)  }\right) \\
&  \ \ \ \ \ \ \ \ \ \ +\sum_{\substack{f:\left[  m\right]  \rightarrow\left[
n\right]  ;\\\text{not }f\left(  \left[  m\right]  \right)  \subseteq
I}}\underbrace{\left[  f\left(  \left[  m\right]  \right)  \subseteq I\right]
}_{\substack{=0\\\text{(since we don't have }f\left(  \left[  m\right]
\right)  \subseteq I\text{)}}}\left(  v_{f\left(  1\right)  }v_{f\left(
2\right)  }\cdots v_{f\left(  m\right)  }\right) \\
&  \ \ \ \ \ \ \ \ \ \ \left(
\begin{array}
[c]{c}%
\text{since each map }f:\left[  m\right]  \rightarrow\left[  n\right]  \text{
satisfies either }f\left(  \left[  m\right]  \right)  \subseteq I\\
\text{or }\left(  \text{not }f\left(  \left[  m\right]  \right)  \subseteq
I\right)  \text{ (but not both)}%
\end{array}
\right) \\
&  =\sum_{\substack{f:\left[  m\right]  \rightarrow\left[  n\right]
;\\f\left(  \left[  m\right]  \right)  \subseteq I}}v_{f\left(  1\right)
}v_{f\left(  2\right)  }\cdots v_{f\left(  m\right)  }+\underbrace{\sum
_{\substack{f:\left[  m\right]  \rightarrow\left[  n\right]  ;\\\text{not
}f\left(  \left[  m\right]  \right)  \subseteq I}}0\left(  v_{f\left(
1\right)  }v_{f\left(  2\right)  }\cdots v_{f\left(  m\right)  }\right)
}_{=0}\\
&  =\sum_{\substack{f:\left[  m\right]  \rightarrow\left[  n\right]
;\\f\left(  \left[  m\right]  \right)  \subseteq I}}v_{f\left(  1\right)
}v_{f\left(  2\right)  }\cdots v_{f\left(  m\right)  }.
\end{align*}
This proves Lemma \ref{lem.sol.noncomm.polarization.1}.
\end{verlong}
\end{proof}

\begin{lemma}
\label{lem.sol.noncomm.polarization.2}Let $n\in\mathbb{N}$. Let $J$ be a
subset of $\left[  n\right]  $. Then,%
\[
\sum_{\substack{I\subseteq\left[  n\right]  ;\\J\subseteq I}}\left(
-1\right)  ^{n-\left\vert I\right\vert }=\left[  J=\left[  n\right]  \right]
.
\]

\end{lemma}

\begin{vershort}
\begin{proof}
[Proof of Lemma \ref{lem.sol.noncomm.polarization.2}.]Here is a proof that
makes use of Exercise \ref{exe.prod(ai+bi)}. It is not the simplest possible
proof, but it might be the shortest.\footnote{We are working in the
commutative ring $\mathbb{Z}$ in this proof. Hence, product signs such as
$\prod_{i\in I}$ and $\prod_{i\in\left[  n\right]  \setminus I}$ make sense.}

We first make the following observations:

\begin{itemize}
\item If $I$ is a subset of $\left[  n\right]  $ satisfying $J\subseteq I$,
then%
\begin{equation}
\prod_{i\in\left[  n\right]  \setminus I}\left[  i\notin J\right]  =1
\label{pf.lem.sol.noncomm.polarization.2.short.1}%
\end{equation}
\footnote{\textit{Proof of (\ref{pf.lem.sol.noncomm.polarization.2.short.1}):}
Let $I$ be a subset of $\left[  n\right]  $ satisfying $J\subseteq I$.
\par
From $J\subseteq I$, we obtain $\left[  n\right]  \setminus I\subseteq\left[
n\right]  \setminus J$. Thus, each $i\in\left[  n\right]  \setminus I$
satisfies $i\in\left[  n\right]  \setminus I\subseteq\left[  n\right]
\setminus J$ and therefore $i\notin J$. Hence, each $i\in\left[  n\right]
\setminus I$ satisfies $\left[  i\notin J\right]  =1$. Therefore, $\prod
_{i\in\left[  n\right]  \setminus I}\underbrace{\left[  i\notin J\right]
}_{=1}=\prod_{i\in\left[  n\right]  \setminus I}1=1$. This proves
(\ref{pf.lem.sol.noncomm.polarization.2.short.1}).}.

\item If $I$ is a subset of $\left[  n\right]  $ that does \textbf{not}
satisfy $J\subseteq I$, then
\begin{equation}
\prod_{i\in\left[  n\right]  \setminus I}\left[  i\notin J\right]  =0
\label{pf.lem.sol.noncomm.polarization.2.short.2}%
\end{equation}
\footnote{\textit{Proof of (\ref{pf.lem.sol.noncomm.polarization.2.short.2}):}
Let $I$ be a subset of $\left[  n\right]  $ that does \textbf{not} satisfy
$J\subseteq I$.
\par
We have $J\not \subseteq I$ (since we do not have $J\subseteq I$). Hence,
there exists some $j\in J$ such that $j\notin I$. Consider this $j$. We have
$\left[  j\notin J\right]  =0$ (since $j\in J$). Combining $j\in
J\subseteq\left[  n\right]  $ with $j\notin I$, we obtain $j\in\left[
n\right]  \setminus I$. Hence, at least one factor of the product $\prod
_{i\in\left[  n\right]  \setminus I}\left[  i\notin J\right]  $ is $0$
(namely, the factor for $i=j$ is $\left[  j\notin J\right]  =0$). Thus, the
whole product $\prod_{i\in\left[  n\right]  \setminus I}\left[  i\notin
J\right]  $ must be $0$ (because if a factor of a product is $0$, then the
whole product must be $0$). This proves
(\ref{pf.lem.sol.noncomm.polarization.2.short.2}).}.
\end{itemize}

But Exercise \ref{exe.prod(ai+bi)} \textbf{(a)} (applied to $\mathbb{Z}$, $1$
and $-\left[  i\notin J\right]  $ instead of $\mathbb{K}$, $a_{i}$ and $b_{i}%
$) yields%
\begin{align}
&  \prod_{i=1}^{n}\left(  1+\left(  -\left[  i\notin J\right]  \right)
\right) \nonumber\\
&  =\sum_{I\subseteq\left[  n\right]  }\underbrace{\left(  \prod_{i\in
I}1\right)  }_{=1}\underbrace{\left(  \prod_{i\in\left[  n\right]  \setminus
I}\left(  -\left[  i\notin J\right]  \right)  \right)  }_{=\left(  -1\right)
^{\left\vert \left[  n\right]  \setminus I\right\vert }\prod_{i\in\left[
n\right]  \setminus I}\left[  i\notin J\right]  }\nonumber\\
&  =\sum_{I\subseteq\left[  n\right]  }\left(  -1\right)  ^{\left\vert \left[
n\right]  \setminus I\right\vert }\prod_{i\in\left[  n\right]  \setminus
I}\left[  i\notin J\right] \nonumber\\
&  =\sum_{\substack{I\subseteq\left[  n\right]  ;\\J\subseteq I}%
}\underbrace{\left(  -1\right)  ^{\left\vert \left[  n\right]  \setminus
I\right\vert }}_{\substack{=\left(  -1\right)  ^{n-\left\vert I\right\vert
}\\\text{(since }\left\vert \left[  n\right]  \setminus I\right\vert
=n-\left\vert I\right\vert \\\text{(since }I\subseteq\left[  n\right]
\text{))}}}\underbrace{\prod_{i\in\left[  n\right]  \setminus I}\left[
i\notin J\right]  }_{\substack{=1\\\text{(by
(\ref{pf.lem.sol.noncomm.polarization.2.short.1}))}}}+\sum
_{\substack{I\subseteq\left[  n\right]  ;\\\text{not }J\subseteq I}}\left(
-1\right)  ^{\left\vert \left[  n\right]  \setminus I\right\vert
}\underbrace{\prod_{i\in\left[  n\right]  \setminus I}\left[  i\notin
J\right]  }_{\substack{=0\\\text{(by
(\ref{pf.lem.sol.noncomm.polarization.2.short.2}))}}}\nonumber\\
&  \ \ \ \ \ \ \ \ \ \ \left(
\begin{array}
[c]{c}%
\text{since each subset }I\text{ of }\left[  n\right]  \text{ satisfies either
}J\subseteq I\text{ or }\left(  \text{not }J\subseteq I\right)  \text{,}\\
\text{but not both}%
\end{array}
\right) \nonumber\\
&  =\sum_{\substack{I\subseteq\left[  n\right]  ;\\J\subseteq I}}\left(
-1\right)  ^{n-\left\vert I\right\vert }+\underbrace{\sum
_{\substack{I\subseteq\left[  n\right]  ;\\\text{not }J\subseteq I}}\left(
-1\right)  ^{\left\vert \left[  n\right]  \setminus I\right\vert }0}_{=0}%
=\sum_{\substack{I\subseteq\left[  n\right]  ;\\J\subseteq I}}\left(
-1\right)  ^{n-\left\vert I\right\vert }.
\label{pf.lem.sol.noncomm.polarization.2.short.8}%
\end{align}

On the other hand, it is easy to see that%
\begin{equation}
\prod_{i=1}^{n}\left(  1+\left(  -\left[  i\notin J\right]  \right)  \right)
=\left[  J=\left[  n\right]  \right]
\label{pf.lem.sol.noncomm.polarization.2.short.9}%
\end{equation}
\footnote{\textit{Proof of (\ref{pf.lem.sol.noncomm.polarization.2.short.9}):}
We are in one of the following two cases:
\par
\textit{Case 1:} We have $J=\left[  n\right]  $.
\par
\textit{Case 2:} We don't have $J=\left[  n\right]  $.
\par
Let us first consider Case 1. In this case, we have $J=\left[  n\right]  $.
Thus, $\left[  J=\left[  n\right]  \right]  =1$. But each $i\in\left\{
1,2,\ldots,n\right\}  $ satisfies $i\in\left\{  1,2,\ldots,n\right\}  =\left[
n\right]  =J$ and therefore $\left[  i\notin J\right]  =0$. Hence,%
\[
\prod_{i=1}^{n}\left(  1+\left(  -\underbrace{\left[  i\notin J\right]  }%
_{=0}\right)  \right)  =\prod_{i=1}^{n}\underbrace{\left(  1+\left(
-0\right)  \right)  }_{=1}=\prod_{i=1}^{n}1=1.
\]
Comparing this with $\left[  J=\left[  n\right]  \right]  =1$, we obtain
$\prod_{i=1}^{n}\left(  1+\left(  -\left[  i\notin J\right]  \right)  \right)
=\left[  J=\left[  n\right]  \right]  $. Hence,
(\ref{pf.lem.sol.noncomm.polarization.2.short.9}) is proven in Case 1.
\par
Let us now consider Case 2. In this case, we don't have $J=\left[  n\right]
$. Hence, $J\neq\left[  n\right]  $. Thus, $J$ is a \textbf{proper} subset of
$\left[  n\right]  $ (since $J\subseteq\left[  n\right]  $). Therefore, there
exists some $j\in\left[  n\right]  $ such that $j\notin J$. Consider this $j$.
We have $\left[  j\notin J\right]  =1$ (since $j\notin J$). But $j\in\left[
n\right]  =\left\{  1,2,\ldots,n\right\}  $. Thus, at least one factor of the
product $\prod_{i=1}^{n}\left(  1+\left(  -\left[  i\notin J\right]  \right)
\right)  $ equals $0$ (namely, the factor for $i=j$ is $1+\left(
-\underbrace{\left[  j\notin J\right]  }_{=1}\right)  =1+\left(  -1\right)
=0$). Consequently, the whole product $\prod_{i=1}^{n}\left(  1+\left(
-\left[  i\notin J\right]  \right)  \right)  $ equals $0$ (because if at least
one factor of a product equals $0$, then the whole product equals $0$). In
other words, $\prod_{i=1}^{n}\left(  1+\left(  -\left[  i\notin J\right]
\right)  \right)  =0$. Comparing this with $\left[  J=\left[  n\right]
\right]  =0$ (which holds because we don't have $J=\left[  n\right]  $), we
obtain $\prod_{i=1}^{n}\left(  1+\left(  -\left[  i\notin J\right]  \right)
\right)  =\left[  J=\left[  n\right]  \right]  $. Hence,
(\ref{pf.lem.sol.noncomm.polarization.2.short.9}) is proven in Case 2.
\par
We have now proven (\ref{pf.lem.sol.noncomm.polarization.2.short.9}) in both
Cases 1 and 2. Thus, (\ref{pf.lem.sol.noncomm.polarization.2.short.9}) is
proven.}. Comparing this with (\ref{pf.lem.sol.noncomm.polarization.2.short.8}%
), we obtain%
\[
\sum_{\substack{I\subseteq\left[  n\right]  ;\\J\subseteq I}}\left(
-1\right)  ^{n-\left\vert I\right\vert }=\left[  J=\left[  n\right]  \right]
.
\]
This proves Lemma \ref{lem.sol.noncomm.polarization.2}.
\end{proof}
\end{vershort}

\begin{verlong}
\begin{proof}
[Proof of Lemma \ref{lem.sol.noncomm.polarization.2}.]There are several ways
to prove Lemma \ref{lem.sol.noncomm.polarization.2} (in particular, there
exist some very simple combinatorial proofs). Let us give one proof (which
might not be the easiest one):\footnote{We are working in the commutative ring
$\mathbb{Z}$ in this proof. Hence, product signs such as $\prod_{i\in I}$ and
$\prod_{i\in\left[  n\right]  \setminus I}$ make sense.}

We first make the following observations:

\begin{itemize}
\item If $I$ is a subset of $\left[  n\right]  $ satisfying $J\subseteq I$,
then%
\begin{equation}
\prod_{i\in\left[  n\right]  \setminus I}\left[  i\notin J\right]  =1
\label{pf.lem.sol.noncomm.polarization.2.1}%
\end{equation}
\footnote{\textit{Proof of (\ref{pf.lem.sol.noncomm.polarization.2.1}):} Let
$I$ be a subset of $\left[  n\right]  $ satisfying $J\subseteq I$.
\par
From $J\subseteq I$, we obtain $\left[  n\right]  \setminus J\supseteq\left[
n\right]  \setminus I$. In other words, $\left[  n\right]  \setminus
I\subseteq\left[  n\right]  \setminus J$. Thus, each $i\in\left[  n\right]
\setminus I$ satisfies $i\in\left[  n\right]  \setminus I\subseteq\left[
n\right]  \setminus J$ and therefore $i\notin J$. Hence, each $i\in\left[
n\right]  \setminus I$ satisfies $\left[  i\notin J\right]  =1$ (since
$i\notin J$). Therefore, $\prod_{i\in\left[  n\right]  \setminus
I}\underbrace{\left[  i\notin J\right]  }_{=1}=\prod_{i\in\left[  n\right]
\setminus I}1=1$. This proves (\ref{pf.lem.sol.noncomm.polarization.2.1}).}.

\item If $I$ is a subset of $\left[  n\right]  $ that does \textbf{not}
satisfy $J\subseteq I$, then
\begin{equation}
\prod_{i\in\left[  n\right]  \setminus I}\left[  i\notin J\right]  =0
\label{pf.lem.sol.noncomm.polarization.2.2}%
\end{equation}
\footnote{\textit{Proof of (\ref{pf.lem.sol.noncomm.polarization.2.2}):} Let
$I$ be a subset of $\left[  n\right]  $ that does \textbf{not} satisfy
$J\subseteq I$.
\par
We have $J\not \subseteq I$ (since we do not have $J\subseteq I$). Hence,
there exists some $j\in J$ such that $j\notin I$. Consider this $j$. We have
$\left[  j\notin J\right]  =0$ (since $j\notin J$ is false (since $j\in J$)).
Combining $j\in J\subseteq\left[  n\right]  $ with $j\notin I$, we obtain
$j\in\left[  n\right]  \setminus I$. Hence, $\left[  j\notin J\right]  $ is a
factor of the product $\prod_{i\in\left[  n\right]  \setminus I}\left[
i\notin J\right]  $ (namely, the factor for $i=j$). This factor is $0$ (since
$\left[  j\notin J\right]  =0$). Therefore, at least one factor of the product
$\prod_{i\in\left[  n\right]  \setminus I}\left[  i\notin J\right]  $ is $0$
(namely, the factor $\left[  j\notin J\right]  $). Thus, the whole product
$\prod_{i\in\left[  n\right]  \setminus I}\left[  i\notin J\right]  $ must be
$0$ (because if a factor of a product is $0$, then the whole product must be
$0$). In other words, $\prod_{i\in\left[  n\right]  \setminus I}\left[
i\notin J\right]  =0$. This proves (\ref{pf.lem.sol.noncomm.polarization.2.2}%
).}.

\item If $I$ is a subset of $\left[  n\right]  $, then%
\begin{align*}
\left\vert \left[  n\right]  \setminus I\right\vert  &
=\underbrace{\left\vert \left[  n\right]  \right\vert }_{=n}-\left\vert
I\right\vert \ \ \ \ \ \ \ \ \ \ \left(  \text{since }I\subseteq\left[
n\right]  \right) \\
&  =n-\left\vert I\right\vert
\end{align*}
and%
\begin{align}
\prod_{i\in\left[  n\right]  \setminus I}\left(  -\left[  i\notin J\right]
\right)   &  =\underbrace{\left(  -1\right)  ^{\left\vert \left[  n\right]
\setminus I\right\vert }}_{\substack{=\left(  -1\right)  ^{n-\left\vert
I\right\vert }\\\text{(since }\left\vert \left[  n\right]  \setminus
I\right\vert =n-\left\vert I\right\vert \text{)}}}\prod_{i\in\left[  n\right]
\setminus I}\left[  i\notin J\right] \nonumber\\
&  =\left(  -1\right)  ^{n-\left\vert I\right\vert }\prod_{i\in\left[
n\right]  \setminus I}\left[  i\notin J\right]  .
\label{pf.lem.sol.noncomm.polarization.2.3}%
\end{align}

\end{itemize}

But Exercise \ref{exe.prod(ai+bi)} \textbf{(a)} (applied to $\mathbb{Z}$, $1$
and $-\left[  i\notin J\right]  $ instead of $\mathbb{K}$, $a_{i}$ and $b_{i}%
$) yields%
\begin{align}
&  \prod_{i=1}^{n}\left(  1+\left(  -\left[  i\notin J\right]  \right)
\right) \nonumber\\
&  =\sum_{I\subseteq\left[  n\right]  }\underbrace{\left(  \prod_{i\in
I}1\right)  }_{=1}\underbrace{\left(  \prod_{i\in\left[  n\right]  \setminus
I}\left(  -\left[  i\notin J\right]  \right)  \right)  }_{\substack{=\left(
-1\right)  ^{n-\left\vert I\right\vert }\prod_{i\in\left[  n\right]  \setminus
I}\left[  i\notin J\right]  \\\text{(by
(\ref{pf.lem.sol.noncomm.polarization.2.3}))}}}\nonumber\\
&  =\sum_{I\subseteq\left[  n\right]  }\left(  -1\right)  ^{n-\left\vert
I\right\vert }\prod_{i\in\left[  n\right]  \setminus I}\left[  i\notin
J\right] \nonumber\\
&  =\sum_{\substack{I\subseteq\left[  n\right]  ;\\J\subseteq I}}\left(
-1\right)  ^{n-\left\vert I\right\vert }\underbrace{\prod_{i\in\left[
n\right]  \setminus I}\left[  i\notin J\right]  }_{\substack{=1\\\text{(by
(\ref{pf.lem.sol.noncomm.polarization.2.1}))}}}+\sum_{\substack{I\subseteq
\left[  n\right]  ;\\\text{not }J\subseteq I}}\left(  -1\right)
^{n-\left\vert I\right\vert }\underbrace{\prod_{i\in\left[  n\right]
\setminus I}\left[  i\notin J\right]  }_{\substack{=0\\\text{(by
(\ref{pf.lem.sol.noncomm.polarization.2.2}))}}}\nonumber\\
&  \ \ \ \ \ \ \ \ \ \ \left(
\begin{array}
[c]{c}%
\text{since each subset }I\text{ of }\left[  n\right]  \text{ satisfies either
}J\subseteq I\text{ or }\left(  \text{not }J\subseteq I\right)  \text{,}\\
\text{but not both}%
\end{array}
\right) \nonumber\\
&  =\sum_{\substack{I\subseteq\left[  n\right]  ;\\J\subseteq I}}\left(
-1\right)  ^{n-\left\vert I\right\vert }1+\underbrace{\sum
_{\substack{I\subseteq\left[  n\right]  ;\\\text{not }J\subseteq I}}\left(
-1\right)  ^{n-\left\vert I\right\vert }0}_{=0}=\sum_{\substack{I\subseteq
\left[  n\right]  ;\\J\subseteq I}}\left(  -1\right)  ^{n-\left\vert
I\right\vert }1\nonumber\\
&  =\sum_{\substack{I\subseteq\left[  n\right]  ;\\J\subseteq I}}\left(
-1\right)  ^{n-\left\vert I\right\vert }.
\label{pf.lem.sol.noncomm.polarization.2.8}%
\end{align}

On the other hand, it is easy to see that%
\begin{equation}
\prod_{i=1}^{n}\left(  1+\left(  -\left[  i\notin J\right]  \right)  \right)
=\left[  J=\left[  n\right]  \right]
\label{pf.lem.sol.noncomm.polarization.2.9}%
\end{equation}
\footnote{\textit{Proof of (\ref{pf.lem.sol.noncomm.polarization.2.9}):} We
are in one of the following two cases:
\par
\textit{Case 1:} We have $J=\left[  n\right]  $.
\par
\textit{Case 2:} We don't have $J=\left[  n\right]  $.
\par
Let us first consider Case 1. In this case, we have $J=\left[  n\right]  $.
Thus, $\left[  J=\left[  n\right]  \right]  =1$. But each $i\in\left\{
1,2,\ldots,n\right\}  $ satisfies $i\in\left\{  1,2,\ldots,n\right\}  =\left[
n\right]  =J$ and therefore $\left[  i\notin J\right]  =0$ (since $i\notin J$
is false (since we have $i\in J$)). Hence,%
\[
\prod_{i=1}^{n}\left(  1+\left(  -\underbrace{\left[  i\notin J\right]  }%
_{=0}\right)  \right)  =\prod_{i=1}^{n}\underbrace{\left(  1+\left(
-0\right)  \right)  }_{=1}=\prod_{i=1}^{n}1=1.
\]
Comparing this with $\left[  J=\left[  n\right]  \right]  =1$, we obtain
$\prod_{i=1}^{n}\left(  1+\left(  -\left[  i\notin J\right]  \right)  \right)
=\left[  J=\left[  n\right]  \right]  $. Hence,
(\ref{pf.lem.sol.noncomm.polarization.2.9}) is proven in Case 1.
\par
Let us now consider Case 2. In this case, we don't have $J=\left[  n\right]
$. Hence, $J\neq\left[  n\right]  $. Thus, $J$ is a \textbf{proper} subset of
$\left[  n\right]  $ (since $J\subseteq\left[  n\right]  $). Therefore, there
exists some $j\in\left[  n\right]  $ such that $j\notin J$. Consider this $j$.
We have $\left[  j\notin J\right]  =1$ (since $j\notin J$). But $j\in\left[
n\right]  =\left\{  1,2,\ldots,n\right\}  $. Thus, $1+\left(  -\left[  j\notin
J\right]  \right)  $ is a factor of the product $\prod_{i=1}^{n}\left(
1+\left(  -\left[  i\notin J\right]  \right)  \right)  $ (namely, the factor
for $i=j$). Moreover, this factor equals $0$ (since $1+\left(
-\underbrace{\left[  j\notin J\right]  }_{=1}\right)  =1+\left(  -1\right)
=0$). Hence, at least one factor of the product $\prod_{i=1}^{n}\left(
1+\left(  -\left[  i\notin J\right]  \right)  \right)  $ equals $0$ (namely,
the factor $1+\left(  -\left[  j\notin J\right]  \right)  $). Consequently,
the whole product $\prod_{i=1}^{n}\left(  1+\left(  -\left[  i\notin J\right]
\right)  \right)  $ equals $0$ (because if at least one factor of a product
equals $0$, then the whole product equals $0$). In other words, $\prod
_{i=1}^{n}\left(  1+\left(  -\left[  i\notin J\right]  \right)  \right)  =0$.
Comparing this with $\left[  J=\left[  n\right]  \right]  =0$ (which holds
because we don't have $J=\left[  n\right]  $), we obtain $\prod_{i=1}%
^{n}\left(  1+\left(  -\left[  i\notin J\right]  \right)  \right)  =\left[
J=\left[  n\right]  \right]  $. Hence,
(\ref{pf.lem.sol.noncomm.polarization.2.9}) is proven in Case 2.
\par
We have now proven (\ref{pf.lem.sol.noncomm.polarization.2.9}) in both Cases 1
and 2. Since these are the only possible cases, we thus have always proven
(\ref{pf.lem.sol.noncomm.polarization.2.9}).}. Comparing this with
(\ref{pf.lem.sol.noncomm.polarization.2.8}), we obtain%
\[
\sum_{\substack{I\subseteq\left[  n\right]  ;\\J\subseteq I}}\left(
-1\right)  ^{n-\left\vert I\right\vert }=\left[  J=\left[  n\right]  \right]
.
\]
This proves Lemma \ref{lem.sol.noncomm.polarization.2}.
\end{proof}
\end{verlong}

\begin{proof}
[Solution to Exercise \ref{exe.noncomm.polarization}.]\textbf{(a)} Let
$m\in\mathbb{N}$.

\begin{vershort}
We have%
\begin{align*}
&  \sum_{I\subseteq\left[  n\right]  }\left(  -1\right)  ^{n-\left\vert
I\right\vert }\underbrace{\left(  \sum_{i\in I}v_{i}\right)  ^{m}%
}_{\substack{=\sum_{\substack{f:\left[  m\right]  \rightarrow\left[  n\right]
;\\f\left(  \left[  m\right]  \right)  \subseteq I}}v_{f\left(  1\right)
}v_{f\left(  2\right)  }\cdots v_{f\left(  m\right)  }\\\text{(by Lemma
\ref{lem.sol.noncomm.polarization.1})}}}\\
&  =\sum_{I\subseteq\left[  n\right]  }\left(  -1\right)  ^{n-\left\vert
I\right\vert }\sum_{\substack{f:\left[  m\right]  \rightarrow\left[  n\right]
;\\f\left(  \left[  m\right]  \right)  \subseteq I}}v_{f\left(  1\right)
}v_{f\left(  2\right)  }\cdots v_{f\left(  m\right)  }\\
&  =\underbrace{\sum_{I\subseteq\left[  n\right]  }\sum_{\substack{f:\left[
m\right]  \rightarrow\left[  n\right]  ;\\f\left(  \left[  m\right]  \right)
\subseteq I}}}_{=\sum_{f:\left[  m\right]  \rightarrow\left[  n\right]  }%
\sum_{\substack{I\subseteq\left[  n\right]  ;\\f\left(  \left[  m\right]
\right)  \subseteq I}}}\left(  -1\right)  ^{n-\left\vert I\right\vert
}v_{f\left(  1\right)  }v_{f\left(  2\right)  }\cdots v_{f\left(  m\right)
}\\
&  =\sum_{f:\left[  m\right]  \rightarrow\left[  n\right]  }\sum
_{\substack{I\subseteq\left[  n\right]  ;\\f\left(  \left[  m\right]  \right)
\subseteq I}}\left(  -1\right)  ^{n-\left\vert I\right\vert }v_{f\left(
1\right)  }v_{f\left(  2\right)  }\cdots v_{f\left(  m\right)  }\\
&  =\sum_{f:\left[  m\right]  \rightarrow\left[  n\right]  }%
\underbrace{\left(  \sum_{\substack{I\subseteq\left[  n\right]  ;\\f\left(
\left[  m\right]  \right)  \subseteq I}}\left(  -1\right)  ^{n-\left\vert
I\right\vert }\right)  }_{\substack{=\left[  f\left(  \left[  m\right]
\right)  =\left[  n\right]  \right]  \\\text{(by Lemma
\ref{lem.sol.noncomm.polarization.2}}\\\text{(applied to }J=f\left(  \left[
m\right]  \right)  \text{))}}}v_{f\left(  1\right)  }v_{f\left(  2\right)
}\cdots v_{f\left(  m\right)  }%
\end{align*}%
\begin{align*}
&  =\sum_{f:\left[  m\right]  \rightarrow\left[  n\right]  }\left[
\underbrace{f\left(  \left[  m\right]  \right)  =\left[  n\right]
}_{\Longleftrightarrow\ \left(  f\text{ is surjective}\right)  }\right]
v_{f\left(  1\right)  }v_{f\left(  2\right)  }\cdots v_{f\left(  m\right)  }\\
&  =\sum_{f:\left[  m\right]  \rightarrow\left[  n\right]  }\left[  f\text{ is
surjective}\right]  v_{f\left(  1\right)  }v_{f\left(  2\right)  }\cdots
v_{f\left(  m\right)  }\\
&  =\sum_{\substack{f:\left[  m\right]  \rightarrow\left[  n\right]
;\\f\text{ is surjective}}}\underbrace{\left[  f\text{ is surjective}\right]
}_{\substack{=1\\\text{(since }f\text{ is surjective)}}}v_{f\left(  1\right)
}v_{f\left(  2\right)  }\cdots v_{f\left(  m\right)  }\\
&  \ \ \ \ \ \ \ \ \ \ +\sum_{\substack{f:\left[  m\right]  \rightarrow\left[
n\right]  ;\\f\text{ is not surjective}}}\underbrace{\left[  f\text{ is
surjective}\right]  }_{\substack{=0\\\text{(since }f\text{ is not
surjective)}}}v_{f\left(  1\right)  }v_{f\left(  2\right)  }\cdots v_{f\left(
m\right)  }\\
&  =\sum_{\substack{f:\left[  m\right]  \rightarrow\left[  n\right]
;\\f\text{ is surjective}}}v_{f\left(  1\right)  }v_{f\left(  2\right)
}\cdots v_{f\left(  m\right)  }+\underbrace{\sum_{\substack{f:\left[
m\right]  \rightarrow\left[  n\right]  ;\\f\text{ is not surjective}%
}}0v_{f\left(  1\right)  }v_{f\left(  2\right)  }\cdots v_{f\left(  m\right)
}}_{=0}\\
&  =\sum_{\substack{f:\left[  m\right]  \rightarrow\left[  n\right]
;\\f\text{ is surjective}}}v_{f\left(  1\right)  }v_{f\left(  2\right)
}\cdots v_{f\left(  m\right)  }.
\end{align*}
This solves Exercise \ref{exe.noncomm.polarization} \textbf{(a)}.
\end{vershort}

\begin{verlong}
We have%
\begin{align*}
&  \sum_{I\subseteq\left[  n\right]  }\left(  -1\right)  ^{n-\left\vert
I\right\vert }\underbrace{\left(  \sum_{i\in I}v_{i}\right)  ^{m}%
}_{\substack{=\sum_{\substack{f:\left[  m\right]  \rightarrow\left[  n\right]
;\\f\left(  \left[  m\right]  \right)  \subseteq I}}v_{f\left(  1\right)
}v_{f\left(  2\right)  }\cdots v_{f\left(  m\right)  }\\\text{(by Lemma
\ref{lem.sol.noncomm.polarization.1})}}}\\
&  =\sum_{I\subseteq\left[  n\right]  }\left(  -1\right)  ^{n-\left\vert
I\right\vert }\sum_{\substack{f:\left[  m\right]  \rightarrow\left[  n\right]
;\\f\left(  \left[  m\right]  \right)  \subseteq I}}v_{f\left(  1\right)
}v_{f\left(  2\right)  }\cdots v_{f\left(  m\right)  }\\
&  =\underbrace{\sum_{I\subseteq\left[  n\right]  }\sum_{\substack{f:\left[
m\right]  \rightarrow\left[  n\right]  ;\\f\left(  \left[  m\right]  \right)
\subseteq I}}}_{=\sum_{f:\left[  m\right]  \rightarrow\left[  n\right]  }%
\sum_{\substack{I\subseteq\left[  n\right]  ;\\f\left(  \left[  m\right]
\right)  \subseteq I}}}\left(  -1\right)  ^{n-\left\vert I\right\vert
}v_{f\left(  1\right)  }v_{f\left(  2\right)  }\cdots v_{f\left(  m\right)
}\\
&  =\sum_{f:\left[  m\right]  \rightarrow\left[  n\right]  }\sum
_{\substack{I\subseteq\left[  n\right]  ;\\f\left(  \left[  m\right]  \right)
\subseteq I}}\left(  -1\right)  ^{n-\left\vert I\right\vert }v_{f\left(
1\right)  }v_{f\left(  2\right)  }\cdots v_{f\left(  m\right)  }\\
&  =\sum_{f:\left[  m\right]  \rightarrow\left[  n\right]  }%
\underbrace{\left(  \sum_{\substack{I\subseteq\left[  n\right]  ;\\f\left(
\left[  m\right]  \right)  \subseteq I}}\left(  -1\right)  ^{n-\left\vert
I\right\vert }\right)  }_{\substack{=\left[  f\left(  \left[  m\right]
\right)  =\left[  n\right]  \right]  \\\text{(by Lemma
\ref{lem.sol.noncomm.polarization.2}}\\\text{(applied to }J=f\left(  \left[
m\right]  \right)  \text{))}}}v_{f\left(  1\right)  }v_{f\left(  2\right)
}\cdots v_{f\left(  m\right)  }%
\end{align*}%
\begin{align*}
&  =\sum_{f:\left[  m\right]  \rightarrow\left[  n\right]  }\left[
\underbrace{f\left(  \left[  m\right]  \right)  =\left[  n\right]
}_{\substack{\Longleftrightarrow\ \left(  f\text{ is surjective}\right)
\\\text{(because }f\text{ is surjective if and only if }f\left(  \left[
m\right]  \right)  =\left[  n\right]  \text{)}}}\right]  v_{f\left(  1\right)
}v_{f\left(  2\right)  }\cdots v_{f\left(  m\right)  }\\
&  =\sum_{f:\left[  m\right]  \rightarrow\left[  n\right]  }\left[  f\text{ is
surjective}\right]  v_{f\left(  1\right)  }v_{f\left(  2\right)  }\cdots
v_{f\left(  m\right)  }\\
&  =\sum_{\substack{f:\left[  m\right]  \rightarrow\left[  n\right]
;\\f\text{ is surjective}}}\underbrace{\left[  f\text{ is surjective}\right]
}_{\substack{=1\\\text{(since }f\text{ is surjective)}}}v_{f\left(  1\right)
}v_{f\left(  2\right)  }\cdots v_{f\left(  m\right)  }\\
&  \ \ \ \ \ \ \ \ \ \ +\sum_{\substack{f:\left[  m\right]  \rightarrow\left[
n\right]  ;\\f\text{ is not surjective}}}\underbrace{\left[  f\text{ is
surjective}\right]  }_{\substack{=0\\\text{(since }f\text{ is not
surjective)}}}v_{f\left(  1\right)  }v_{f\left(  2\right)  }\cdots v_{f\left(
m\right)  }\\
&  \ \ \ \ \ \ \ \ \ \ \left(
\begin{array}
[c]{c}%
\text{because for each }f:\left[  m\right]  \rightarrow\left[  n\right]
\text{, either }\left(  f\text{ is surjective}\right) \\
\text{or }\left(  f\text{ is not surjective}\right)  \text{ (but not both)}%
\end{array}
\right) \\
&  =\sum_{\substack{f:\left[  m\right]  \rightarrow\left[  n\right]
;\\f\text{ is surjective}}}v_{f\left(  1\right)  }v_{f\left(  2\right)
}\cdots v_{f\left(  m\right)  }+\underbrace{\sum_{\substack{f:\left[
m\right]  \rightarrow\left[  n\right]  ;\\f\text{ is not surjective}%
}}0v_{f\left(  1\right)  }v_{f\left(  2\right)  }\cdots v_{f\left(  m\right)
}}_{=0}\\
&  =\sum_{\substack{f:\left[  m\right]  \rightarrow\left[  n\right]
;\\f\text{ is surjective}}}v_{f\left(  1\right)  }v_{f\left(  2\right)
}\cdots v_{f\left(  m\right)  }.
\end{align*}
This solves Exercise \ref{exe.noncomm.polarization} \textbf{(a)}.
\end{verlong}

\textbf{(b)} Let $m\in\left\{  0,1,\ldots,n-1\right\}  $. Then, there exists
no map $f:\left[  m\right]  \rightarrow\left[  n\right]  $ such that $f$ is
surjective\footnote{\textit{Proof.} Let $f:\left[  m\right]  \rightarrow
\left[  n\right]  $ be a map such that $f$ is surjective. Thus, $f\left(
\left[  m\right]  \right)  =\left[  n\right]  $ (since $f$ is surjective). But
clearly, $\left\vert f\left(  \left[  m\right]  \right)  \right\vert
\leq\left\vert \left[  m\right]  \right\vert =m$. Since $f\left(  \left[
m\right]  \right)  =\left[  n\right]  $, this rewrites as $\left\vert \left[
n\right]  \right\vert \leq m$. Since $\left\vert \left[  n\right]  \right\vert
=n$, this rewrites as $n\leq m$. But $m\leq n-1$ (since $m\in\left\{
0,1,\ldots,n-1\right\}  $). Thus, $n\leq m\leq n-1<n$. This is absurd.
\par
Now, forget that we fixed $f$. We thus have found a contradiction for each map
$f:\left[  m\right]  \rightarrow\left[  n\right]  $ such that $f$ is
surjective. Hence, there exists no map $f:\left[  m\right]  \rightarrow\left[
n\right]  $ such that $f$ is surjective. Qed.}. Hence, the sum $\sum
_{\substack{f:\left[  m\right]  \rightarrow\left[  n\right]  ;\\f\text{ is
surjective}}}v_{f\left(  1\right)  }v_{f\left(  2\right)  }\cdots v_{f\left(
m\right)  }$ is an empty sum. Therefore, this sum equals $0$. In other words,
$\sum_{\substack{f:\left[  m\right]  \rightarrow\left[  n\right]  ;\\f\text{
is surjective}}}v_{f\left(  1\right)  }v_{f\left(  2\right)  }\cdots
v_{f\left(  m\right)  }=0$.

Now, Exercise \ref{exe.noncomm.polarization} \textbf{(a)} yields%
\[
\sum_{I\subseteq\left[  n\right]  }\left(  -1\right)  ^{n-\left\vert
I\right\vert }\left(  \sum_{i\in I}v_{i}\right)  ^{m}=\sum
_{\substack{f:\left[  m\right]  \rightarrow\left[  n\right]  ;\\f\text{ is
surjective}}}v_{f\left(  1\right)  }v_{f\left(  2\right)  }\cdots v_{f\left(
m\right)  }=0.
\]
This solves Exercise \ref{exe.noncomm.polarization} \textbf{(b)}.

\textbf{(c)} Clearly, $\left[  n\right]  $ and $\left[  n\right]  $ are two
finite sets such that $\left\vert \left[  n\right]  \right\vert \leq\left\vert
\left[  n\right]  \right\vert $. Hence, Lemma \ref{lem.jectivity.pigeon-surj}
(applied to $U=\left[  n\right]  $ and $V=\left[  n\right]  $) shows that if
$f:\left[  n\right]  \rightarrow\left[  n\right]  $ is a map, then we have the
following logical equivalence:%
\begin{equation}
\left(  f\text{ is surjective}\right)  \ \Longleftrightarrow\ \left(  f\text{
is bijective}\right)  . \label{sol.noncomm.polarization.c.3}%
\end{equation}

\begin{vershort}
Hence, we have the following equality of summation signs:%
\[
\sum_{\substack{f:\left[  n\right]  \rightarrow\left[  n\right]  ;\\f\text{ is
surjective}}}=\sum_{\substack{f:\left[  n\right]  \rightarrow\left[  n\right]
;\\f\text{ is bijective}}}=\sum_{f\text{ is a permutation of }\left[
n\right]  }=\sum_{f\in S_{n}}%
\]
(since $S_{n}$ is the set of all permutations of $\left[  n\right]  $). Now,
Exercise \ref{exe.noncomm.polarization} \textbf{(a)} (applied to $m=n$) yields%
\begin{align*}
\sum_{I\subseteq\left[  n\right]  }\left(  -1\right)  ^{n-\left\vert
I\right\vert }\left(  \sum_{i\in I}v_{i}\right)  ^{n}  &  =\underbrace{\sum
_{\substack{f:\left[  n\right]  \rightarrow\left[  n\right]  ;\\f\text{ is
surjective}}}}_{\substack{=\sum_{f\in S_{n}}\\\text{(by
(\ref{sol.noncomm.polarization.c.3}))}}}v_{f\left(  1\right)  }v_{f\left(
2\right)  }\cdots v_{f\left(  n\right)  }\\
&  =\sum_{f\in S_{n}}v_{f\left(  1\right)  }v_{f\left(  2\right)  }\cdots
v_{f\left(  n\right)  }=\sum_{\sigma\in S_{n}}v_{\sigma\left(  1\right)
}v_{\sigma\left(  2\right)  }\cdots v_{\sigma\left(  n\right)  }%
\end{align*}
(here, we have renamed the summation index $f$ as $\sigma$). This solves
Exercise \ref{exe.noncomm.polarization} \textbf{(c)}.
\end{vershort}

\begin{verlong}
Now, Exercise \ref{exe.noncomm.polarization} \textbf{(a)} (applied to $m=n$)
yields%
\begin{align}
\sum_{I\subseteq\left[  n\right]  }\left(  -1\right)  ^{n-\left\vert
I\right\vert }\left(  \sum_{i\in I}v_{i}\right)  ^{n}  &  =\underbrace{\sum
_{\substack{f:\left[  n\right]  \rightarrow\left[  n\right]  ;\\f\text{ is
surjective}}}}_{\substack{=\sum_{\substack{f:\left[  n\right]  \rightarrow
\left[  n\right]  ;\\f\text{ is bijective}}}\\\text{(because for a map
}f:\left[  n\right]  \rightarrow\left[  n\right]  \text{,}\\\text{the
condition }\left(  f\text{ is surjective}\right)  \text{ is equivalent}%
\\\text{to the condition }\left(  f\text{ is bijective}\right)
\\\text{(because of (\ref{sol.noncomm.polarization.c.3})))}}}v_{f\left(
1\right)  }v_{f\left(  2\right)  }\cdots v_{f\left(  n\right)  }\nonumber\\
&  =\sum_{\substack{f:\left[  n\right]  \rightarrow\left[  n\right]
;\\f\text{ is bijective}}}v_{f\left(  1\right)  }v_{f\left(  2\right)  }\cdots
v_{f\left(  n\right)  }. \label{sol.noncomm.polarization.c.5}%
\end{align}

But $S_{n}$ is the set of all permutations of the set $\left\{  1,2,\ldots
,n\right\}  $. In other words, $S_{n}$ is the set of all permutations of the
set $\left[  n\right]  $ (since $\left\{  1,2,\ldots,n\right\}  =\left[
n\right]  $). In other words, $S_{n}$ is the set of all bijective maps
$\left[  n\right]  \rightarrow\left[  n\right]  $ (since the permutations of
the set $\left[  n\right]  $ are exactly the bijective maps $\left[  n\right]
\rightarrow\left[  n\right]  $). Hence,%
\begin{equation}
\sum_{f\in S_{n}}=\sum_{f\text{ is a bijective map }\left[  n\right]
\rightarrow\left[  n\right]  }=\sum_{\substack{f:\left[  n\right]
\rightarrow\left[  n\right]  ;\\f\text{ is bijective}}}
\label{sol.noncomm.polarization.c.7}%
\end{equation}
(an equality of summation signs). Now, (\ref{sol.noncomm.polarization.c.5})
becomes%
\begin{align*}
\sum_{I\subseteq\left[  n\right]  }\left(  -1\right)  ^{n-\left\vert
I\right\vert }\left(  \sum_{i\in I}v_{i}\right)  ^{n}  &  =\underbrace{\sum
_{\substack{f:\left[  n\right]  \rightarrow\left[  n\right]  ;\\f\text{ is
bijective}}}}_{\substack{=\sum_{f\in S_{n}}\\\text{(by
(\ref{sol.noncomm.polarization.c.7}))}}}v_{f\left(  1\right)  }v_{f\left(
2\right)  }\cdots v_{f\left(  n\right)  }=\sum_{f\in S_{n}}v_{f\left(
1\right)  }v_{f\left(  2\right)  }\cdots v_{f\left(  n\right)  }\\
&  =\sum_{\sigma\in S_{n}}v_{\sigma\left(  1\right)  }v_{\sigma\left(
2\right)  }\cdots v_{\sigma\left(  n\right)  }%
\end{align*}
(here, we have renamed the summation index $f$ as $\sigma$). This solves
Exercise \ref{exe.noncomm.polarization} \textbf{(c)}.
\end{verlong}

\textbf{(d)} We have%
\begin{equation}
v_{\sigma\left(  1\right)  }v_{\sigma\left(  2\right)  }\cdots v_{\sigma
\left(  n\right)  }=v_{1}v_{2}\cdots v_{n}
\label{sol.noncomm.polarization.d.1}%
\end{equation}
for each $\sigma\in S_{n}$\ \ \ \ \footnote{\textit{Proof of
(\ref{sol.noncomm.polarization.d.1}):} Let $\sigma\in S_{n}$. Thus, $\sigma$
is a permutation of the set $\left\{  1,2,\ldots,n\right\}  $ (since $S_{n}$
is the set of all permutations of the set $\left\{  1,2,\ldots,n\right\}  $).
In other words, $\sigma$ is a bijection $\left\{  1,2,\ldots,n\right\}
\rightarrow\left\{  1,2,\ldots,n\right\}  $.
\par
But%
\begin{align*}
v_{\sigma\left(  1\right)  }v_{\sigma\left(  2\right)  }\cdots v_{\sigma
\left(  n\right)  }  &  =\underbrace{\prod_{i=1}^{n}}_{=\prod_{i\in\left\{
1,2,\ldots,n\right\}  }}v_{\sigma\left(  i\right)  }=\prod_{i\in\left\{
1,2,\ldots,n\right\}  }v_{\sigma\left(  i\right)  }=\underbrace{\prod
_{i\in\left\{  1,2,\ldots,n\right\}  }}_{=\prod_{i=1}^{n}}v_{i}\\
&  \ \ \ \ \ \ \ \ \ \ \left(
\begin{array}
[c]{c}%
\text{here, we have substituted }i\text{ for }\sigma\left(  i\right)  \text{
in the product,}\\
\text{since }\sigma\text{ is a bijection }\left\{  1,2,\ldots,n\right\}
\rightarrow\left\{  1,2,\ldots,n\right\}
\end{array}
\right) \\
&  =\prod_{i=1}^{n}v_{i}=v_{1}v_{2}\cdots v_{n}.
\end{align*}
This proves (\ref{sol.noncomm.polarization.d.1}).}. Now, Exercise
\ref{exe.noncomm.polarization} \textbf{(c)} yields%
\begin{align*}
\sum_{I\subseteq\left[  n\right]  }\left(  -1\right)  ^{n-\left\vert
I\right\vert }\left(  \sum_{i\in I}v_{i}\right)  ^{n}  &  =\sum_{\sigma\in
S_{n}}\underbrace{v_{\sigma\left(  1\right)  }v_{\sigma\left(  2\right)
}\cdots v_{\sigma\left(  n\right)  }}_{\substack{=v_{1}v_{2}\cdots
v_{n}\\\text{(by (\ref{sol.noncomm.polarization.d.1}))}}}=\sum_{\sigma\in
S_{n}}v_{1}v_{2}\cdots v_{n}\\
&  =\underbrace{\left\vert S_{n}\right\vert }_{=n!}v_{1}v_{2}\cdots
v_{n}=n!\cdot v_{1}v_{2}\cdots v_{n}.
\end{align*}
This solves Exercise \ref{exe.noncomm.polarization} \textbf{(d)}.
\end{proof}

\subsection{\label{sect.sol.noncomm.polarization2}Solution to Exercise
\ref{exe.noncomm.polarization2}}

\begin{vershort}
Throughout Section \ref{sect.sol.noncomm.polarization2}, we shall use the same
notations that we have used in Section \ref{sect.sol.noncomm.polarization}.
Also, for any set $S$, we let $\mathcal{P}\left(  S\right)  $ denote the
powerset of $S$ (that is, the set of all subsets of $S$).
\end{vershort}

\begin{verlong}
Throughout Section \ref{sect.sol.noncomm.polarization2}, we shall use the same
notations that we have used in Section \ref{sect.sol.noncomm.polarization}. We
shall also use Definition \ref{def.sol.prod(ai+bi).powerset}.
\end{verlong}

\begin{proof}
[Solution to Exercise \ref{exe.noncomm.polarization2}.]We have $n+1\in\left[
n+1\right]  $ (since $n+1$ is a positive integer).

\begin{vershort}
Also, $\left[  n+1\right]  \setminus\left\{  n+1\right\}  =\left[  n\right]  $
(for the same reason).
\end{vershort}

\begin{verlong}
Also,%
\begin{align}
\underbrace{\left[  n+1\right]  }_{\substack{=\left\{  1,2,\ldots,n+1\right\}
\\\text{(by the definition of }\left[  n+1\right]  \text{)}}}\setminus\left\{
n+1\right\}   &  =\left\{  1,2,\ldots,n+1\right\}  \setminus\left\{
n+1\right\} \nonumber\\
&  =\left\{  1,2,\ldots,n\right\}  \ \ \ \ \ \ \ \ \ \ \left(  \text{since
}n\in\mathbb{N}\right) \nonumber\\
&  =\left[  n\right]  \label{sol.noncomm.polarization2.1}%
\end{align}
(since $\left[  n\right]  =\left\{  1,2,\ldots,n\right\}  $).
\end{verlong}

\begin{vershort}
But Fact 1 from the solution to Exercise \ref{exe.prod(ai+bi)} (applied to
$S=\left[  n+1\right]  $ and $s=n+1$) yields the following two facts:
\end{vershort}

\begin{verlong}
But Proposition \ref{prop.sol.prod(ai+bi).powerset.lem} (applied to $S=\left[
n+1\right]  $ and $s=n+1$) yields the following two facts:
\end{verlong}

\begin{itemize}
\item We have $\mathcal{P}\left(  \left[  n+1\right]  \setminus\left\{
n+1\right\}  \right)  \subseteq\mathcal{P}\left(  \left[  n+1\right]  \right)
$.

\item The map%
\begin{align*}
\mathcal{P}\left(  \left[  n+1\right]  \setminus\left\{  n+1\right\}  \right)
&  \rightarrow\mathcal{P}\left(  \left[  n+1\right]  \right)  \setminus
\mathcal{P}\left(  \left[  n+1\right]  \setminus\left\{  n+1\right\}  \right)
,\\
U  &  \mapsto U\cup\left\{  n+1\right\}
\end{align*}
is well-defined and a bijection.
\end{itemize}

\begin{vershort}
Using $\left[  n+1\right]  \setminus\left\{  n+1\right\}  =\left[  n\right]
$, we can rewrite these two facts as follows:
\end{vershort}

\begin{verlong}
Using (\ref{sol.noncomm.polarization2.1}), we can rewrite these two facts as follows:
\end{verlong}

\begin{itemize}
\item We have $\mathcal{P}\left(  \left[  n\right]  \right)  \subseteq
\mathcal{P}\left(  \left[  n+1\right]  \right)  $.

\item The map%
\begin{align*}
\mathcal{P}\left(  \left[  n\right]  \right)   &  \rightarrow\mathcal{P}%
\left(  \left[  n+1\right]  \right)  \setminus\mathcal{P}\left(  \left[
n\right]  \right)  ,\\
U  &  \mapsto U\cup\left\{  n+1\right\}
\end{align*}
is well-defined and a bijection.
\end{itemize}

We now claim that%
\begin{equation}
\sum_{I\subseteq\left[  n\right]  }\left(  -1\right)  ^{n-\left\vert
I\right\vert }\left(  w+\sum_{i\in I}v_{i}\right)  ^{m}=\sum_{I\subseteq
\left[  n\right]  }\left(  -1\right)  ^{n-\left\vert I\right\vert }\left(
\sum_{i\in I}v_{i}\right)  ^{m} \label{sol.noncomm.polarization2.wout}%
\end{equation}
for every $m\in\left\{  0,1,\ldots,n\right\}  $ and $w\in\mathbb{L}$.

\begin{vershort}
[\textit{Proof of (\ref{sol.noncomm.polarization2.wout}):} Let $m\in\left\{
0,1,\ldots,n\right\}  $ and $w\in\mathbb{L}$.

We extend the $n$-tuple $\left(  v_{1},v_{2},\ldots,v_{n}\right)
\in\mathbb{L}^{n}$ to an $\left(  n+1\right)  $-tuple $\left(  v_{1}%
,v_{2},\ldots,v_{n+1}\right)  \in\mathbb{L}^{n+1}$ by setting $v_{n+1}=w$.
Thus, $v_{1},v_{2},\ldots,v_{n+1}$ are $n+1$ elements of $\mathbb{L}$.
Moreover, $m\in\left\{  0,1,\ldots,n\right\}  =\left\{  0,1,\ldots,\left(
n+1\right)  -1\right\}  $. Hence, Exercise \ref{exe.noncomm.polarization}
\textbf{(b)} (applied to $n+1$ instead of $n$) yields%
\[
\sum_{I\subseteq\left[  n+1\right]  }\left(  -1\right)  ^{\left(  n+1\right)
-\left\vert I\right\vert }\left(  \sum_{i\in I}v_{i}\right)  ^{m}=0.
\]
Hence,%
\begin{align*}
0  &  =\underbrace{\sum_{I\subseteq\left[  n+1\right]  }}_{\substack{=\sum
_{I\in\mathcal{P}\left(  \left[  n+1\right]  \right)  }}}\left(  -1\right)
^{\left(  n+1\right)  -\left\vert I\right\vert }\left(  \sum_{i\in I}%
v_{i}\right)  ^{m}=\sum_{I\in\mathcal{P}\left(  \left[  n+1\right]  \right)
}\left(  -1\right)  ^{\left(  n+1\right)  -\left\vert I\right\vert }\left(
\sum_{i\in I}v_{i}\right)  ^{m}\\
&  =\underbrace{\sum_{\substack{I\in\mathcal{P}\left(  \left[  n+1\right]
\right)  ;\\I\in\mathcal{P}\left(  \left[  n\right]  \right)  }}}%
_{\substack{=\sum_{I\in\mathcal{P}\left(  \left[  n\right]  \right)
}\\\text{(since }\mathcal{P}\left(  \left[  n\right]  \right)  \subseteq
\mathcal{P}\left(  \left[  n+1\right]  \right)  \text{)}}}\underbrace{\left(
-1\right)  ^{\left(  n+1\right)  -\left\vert I\right\vert }}%
_{\substack{=-\left(  -1\right)  ^{n-\left\vert I\right\vert }\\\text{(since
}\left(  n+1\right)  -\left\vert I\right\vert =\left(  n-\left\vert
I\right\vert \right)  +1\text{)}}}\left(  \sum_{i\in I}v_{i}\right)  ^{m}\\
&  \ \ \ \ \ \ \ \ \ \ +\underbrace{\sum_{\substack{I\in\mathcal{P}\left(
\left[  n+1\right]  \right)  ;\\I\notin\mathcal{P}\left(  \left[  n\right]
\right)  }}}_{=\sum_{I\in\mathcal{P}\left(  \left[  n+1\right]  \right)
\setminus\mathcal{P}\left(  \left[  n\right]  \right)  }}\left(  -1\right)
^{\left(  n+1\right)  -\left\vert I\right\vert }\left(  \sum_{i\in I}%
v_{i}\right)  ^{m}\\
&  \ \ \ \ \ \ \ \ \ \ \left(
\begin{array}
[c]{c}%
\text{since each }I\in\mathcal{P}\left(  \left[  n+1\right]  \right)  \text{
satisfies either }I\in\mathcal{P}\left(  \left[  n\right]  \right) \\
\text{or }I\notin\mathcal{P}\left(  \left[  n\right]  \right)  \text{ (but not
both)}%
\end{array}
\right) \\
&  =-\sum_{I\in\mathcal{P}\left(  \left[  n\right]  \right)  }\left(
-1\right)  ^{n-\left\vert I\right\vert }\left(  \sum_{i\in I}v_{i}\right)
^{m}+\sum_{I\in\mathcal{P}\left(  \left[  n+1\right]  \right)  \setminus
\mathcal{P}\left(  \left[  n\right]  \right)  }\left(  -1\right)  ^{\left(
n+1\right)  -\left\vert I\right\vert }\left(  \sum_{i\in I}v_{i}\right)  ^{m}.
\end{align*}
Adding $\sum_{I\in\mathcal{P}\left(  \left[  n\right]  \right)  }\left(
-1\right)  ^{n-\left\vert I\right\vert }\left(  \sum_{i\in I}v_{i}\right)
^{m}$ to both sides of this equality, we obtain%
\begin{align}
&  \sum_{I\in\mathcal{P}\left(  \left[  n\right]  \right)  }\left(  -1\right)
^{n-\left\vert I\right\vert }\left(  \sum_{i\in I}v_{i}\right)  ^{m}%
\nonumber\\
&  =\sum_{I\in\mathcal{P}\left(  \left[  n+1\right]  \right)  \setminus
\mathcal{P}\left(  \left[  n\right]  \right)  }\left(  -1\right)  ^{\left(
n+1\right)  -\left\vert I\right\vert }\left(  \sum_{i\in I}v_{i}\right)
^{m}\nonumber\\
&  =\sum_{U\in\mathcal{P}\left(  \left[  n\right]  \right)  }\left(
-1\right)  ^{\left(  n+1\right)  -\left\vert U\cup\left\{  n+1\right\}
\right\vert }\left(  \sum_{i\in U\cup\left\{  n+1\right\}  }v_{i}\right)
^{m}\nonumber\\
&  \ \ \ \ \ \ \ \ \ \ \left(
\begin{array}
[c]{c}%
\text{here, we have substituted }U\cup\left\{  n+1\right\}  \text{ for
}I\text{ in the sum,}\\
\text{since the map }\mathcal{P}\left(  \left[  n\right]  \right)
\rightarrow\mathcal{P}\left(  \left[  n+1\right]  \right)  \setminus
\mathcal{P}\left(  \left[  n\right]  \right)  ,\ U\mapsto U\cup\left\{
n+1\right\} \\
\text{is a bijection}%
\end{array}
\right) \nonumber\\
&  =\sum_{I\in\mathcal{P}\left(  \left[  n\right]  \right)  }\left(
-1\right)  ^{\left(  n+1\right)  -\left\vert I\cup\left\{  n+1\right\}
\right\vert }\left(  \sum_{i\in I\cup\left\{  n+1\right\}  }v_{i}\right)  ^{m}
\label{sol.noncomm.polarization2.wout.pf.short.1}%
\end{align}
(here, we have renamed the summation index $U$ as $I$).

But each $I\in\mathcal{P}\left(  \left[  n\right]  \right)  $ satisfies%
\begin{equation}
\left(  -1\right)  ^{\left(  n+1\right)  -\left\vert I\cup\left\{
n+1\right\}  \right\vert }\left(  \sum_{i\in I\cup\left\{  n+1\right\}  }%
v_{i}\right)  ^{m}=\left(  -1\right)  ^{n-\left\vert I\right\vert }\left(
w+\sum_{i\in I}v_{i}\right)  ^{m}
\label{sol.noncomm.polarization2.wout.pf.short.2}%
\end{equation}
\footnote{\textit{Proof of (\ref{sol.noncomm.polarization2.wout.pf.short.2}):}
Let $I\in\mathcal{P}\left(  \left[  n\right]  \right)  $. Thus, $I$ is a
subset of $\left[  n\right]  $ (since $\mathcal{P}\left(  \left[  n\right]
\right)  $ is the set of all subsets of $\left[  n\right]  $). In other words,
$I\subseteq\left[  n\right]  $. Hence, $n+1\notin I$. Therefore, $\left\vert
I\cup\left\{  n+1\right\}  \right\vert =\left\vert I\right\vert +1$ and
$I\setminus\left\{  n+1\right\}  =I$. Therefore,
\[
\left(  n+1\right)  -\underbrace{\left\vert I\cup\left\{  n+1\right\}
\right\vert }_{=\left\vert I\right\vert +1}=\left(  n+1\right)  -\left(
\left\vert I\right\vert +1\right)  =n-\left\vert I\right\vert .
\]
Furthermore, the sets $I$ and $\left\{  n+1\right\}  $ are disjoint (since
$n+1\notin I$). Thus,%
\[
\sum_{i\in I\cup\left\{  n+1\right\}  }v_{i}=\sum_{i\in I}v_{i}%
+\underbrace{\sum_{i\in\left\{  n+1\right\}  }v_{i}}_{=v_{n+1}=w}=\sum_{i\in
I}v_{i}+w=w+\sum_{i\in I}v_{i}.
\]
Thus,%
\[
\underbrace{\left(  -1\right)  ^{\left(  n+1\right)  -\left\vert I\cup\left\{
n+1\right\}  \right\vert }}_{\substack{=\left(  -1\right)  ^{n-\left\vert
I\right\vert }\\\text{(since }\left(  n+1\right)  -\left\vert I\cup\left\{
n+1\right\}  \right\vert =n-\left\vert I\right\vert \text{)}}}\left(
\underbrace{\sum_{i\in I\cup\left\{  n+1\right\}  }v_{i}}_{=w+\sum_{i\in
I}v_{i}}\right)  ^{m}=\left(  -1\right)  ^{n-\left\vert I\right\vert }\left(
w+\sum_{i\in I}v_{i}\right)  ^{m}.
\]
This proves (\ref{sol.noncomm.polarization2.wout.pf.short.2}).}. Thus,
(\ref{sol.noncomm.polarization2.wout.pf.short.1}) becomes%
\begin{align*}
\sum_{I\in\mathcal{P}\left(  \left[  n\right]  \right)  }\left(  -1\right)
^{n-\left\vert I\right\vert }\left(  \sum_{i\in I}v_{i}\right)  ^{m}  &
=\sum_{I\in\mathcal{P}\left(  \left[  n\right]  \right)  }\underbrace{\left(
-1\right)  ^{\left(  n+1\right)  -\left\vert I\cup\left\{  n+1\right\}
\right\vert }\left(  \sum_{i\in I\cup\left\{  n+1\right\}  }v_{i}\right)
^{m}}_{\substack{=\left(  -1\right)  ^{n-\left\vert I\right\vert }\left(
w+\sum_{i\in I}v_{i}\right)  ^{m}\\\text{(by
(\ref{sol.noncomm.polarization2.wout.pf.short.2}))}}}\\
&  =\sum_{I\in\mathcal{P}\left(  \left[  n\right]  \right)  }\left(
-1\right)  ^{n-\left\vert I\right\vert }\left(  w+\sum_{i\in I}v_{i}\right)
^{m}.
\end{align*}
In view of $\sum_{I\in\mathcal{P}\left(  \left[  n\right]  \right)  }%
=\sum_{I\subseteq\left[  n\right]  }$ (an equality between summation signs),
this equality rewrites as%
\[
\sum_{I\subseteq\left[  n\right]  }\left(  -1\right)  ^{n-\left\vert
I\right\vert }\left(  \sum_{i\in I}v_{i}\right)  ^{m}=\sum_{I\subseteq\left[
n\right]  }\left(  -1\right)  ^{n-\left\vert I\right\vert }\left(
w+\sum_{i\in I}v_{i}\right)  ^{m}.
\]
This proves (\ref{sol.noncomm.polarization2.wout}).]
\end{vershort}

\begin{verlong}
[\textit{Proof of (\ref{sol.noncomm.polarization2.wout}):} Let $m\in\left\{
0,1,\ldots,n\right\}  $ and $w\in\mathbb{L}$.

We extend the $n$-tuple $\left(  v_{1},v_{2},\ldots,v_{n}\right)
\in\mathbb{L}^{n}$ to an $\left(  n+1\right)  $-tuple $\left(  v_{1}%
,v_{2},\ldots,v_{n+1}\right)  \in\mathbb{L}^{n+1}$ by setting $v_{n+1}=w$.
Thus, $v_{1},v_{2},\ldots,v_{n+1}$ are $n+1$ elements of $\mathbb{L}$.
Moreover, $m\in\left\{  0,1,\ldots,n\right\}  =\left\{  0,1,\ldots,\left(
n+1\right)  -1\right\}  $ (since $n=\left(  n+1\right)  -1$). Hence, Exercise
\ref{exe.noncomm.polarization} \textbf{(b)} (applied to $n+1$ instead of $n$)
yields%
\[
\sum_{I\subseteq\left[  n+1\right]  }\left(  -1\right)  ^{\left(  n+1\right)
-\left\vert I\right\vert }\left(  \sum_{i\in I}v_{i}\right)  ^{m}=0.
\]
Hence,%
\begin{align*}
0  &  =\underbrace{\sum_{I\subseteq\left[  n+1\right]  }}_{\substack{=\sum
_{I\in\mathcal{P}\left(  \left[  n+1\right]  \right)  }\\\text{(by the
definition of}\\\text{the summation}\\\text{sign }\sum_{I\subseteq\left[
n+1\right]  }\text{)}}}\left(  -1\right)  ^{\left(  n+1\right)  -\left\vert
I\right\vert }\left(  \sum_{i\in I}v_{i}\right)  ^{m}\\
&  =\sum_{I\in\mathcal{P}\left(  \left[  n+1\right]  \right)  }\left(
-1\right)  ^{\left(  n+1\right)  -\left\vert I\right\vert }\left(  \sum_{i\in
I}v_{i}\right)  ^{m}\\
&  =\underbrace{\sum_{\substack{I\in\mathcal{P}\left(  \left[  n+1\right]
\right)  ;\\I\in\mathcal{P}\left(  \left[  n\right]  \right)  }}}%
_{\substack{=\sum_{I\in\mathcal{P}\left(  \left[  n\right]  \right)
}\\\text{(since }\mathcal{P}\left(  \left[  n\right]  \right)  \subseteq
\mathcal{P}\left(  \left[  n+1\right]  \right)  \text{)}}}\underbrace{\left(
-1\right)  ^{\left(  n+1\right)  -\left\vert I\right\vert }}%
_{\substack{=\left(  -1\right)  ^{\left(  n-\left\vert I\right\vert \right)
+1}\\\text{(since }\left(  n+1\right)  -\left\vert I\right\vert =\left(
n-\left\vert I\right\vert \right)  +1\text{)}}}\left(  \sum_{i\in I}%
v_{i}\right)  ^{m}\\
&  \ \ \ \ \ \ \ \ \ \ +\underbrace{\sum_{\substack{I\in\mathcal{P}\left(
\left[  n+1\right]  \right)  ;\\I\notin\mathcal{P}\left(  \left[  n\right]
\right)  }}}_{=\sum_{I\in\mathcal{P}\left(  \left[  n+1\right]  \right)
\setminus\mathcal{P}\left(  \left[  n\right]  \right)  }}\left(  -1\right)
^{\left(  n+1\right)  -\left\vert I\right\vert }\left(  \sum_{i\in I}%
v_{i}\right)  ^{m}\\
&  \ \ \ \ \ \ \ \ \ \ \left(
\begin{array}
[c]{c}%
\text{since each }I\in\mathcal{P}\left(  \left[  n+1\right]  \right)  \text{
satisfies either }I\in\mathcal{P}\left(  \left[  n\right]  \right) \\
\text{or }I\notin\mathcal{P}\left(  \left[  n\right]  \right)  \text{ (but not
both)}%
\end{array}
\right) \\
&  =\sum_{I\in\mathcal{P}\left(  \left[  n\right]  \right)  }%
\underbrace{\left(  -1\right)  ^{\left(  n-\left\vert I\right\vert \right)
+1}}_{\substack{=\left(  -1\right)  ^{n-\left\vert I\right\vert }\left(
-1\right)  \\=-\left(  -1\right)  ^{n-\left\vert I\right\vert }}}\left(
\sum_{i\in I}v_{i}\right)  ^{m}+\sum_{I\in\mathcal{P}\left(  \left[
n+1\right]  \right)  \setminus\mathcal{P}\left(  \left[  n\right]  \right)
}\left(  -1\right)  ^{\left(  n+1\right)  -\left\vert I\right\vert }\left(
\sum_{i\in I}v_{i}\right)  ^{m}\\
&  =\underbrace{\sum_{I\in\mathcal{P}\left(  \left[  n\right]  \right)
}\left(  -\left(  -1\right)  ^{n-\left\vert I\right\vert }\left(  \sum_{i\in
I}v_{i}\right)  ^{m}\right)  }_{=-\sum_{I\in\mathcal{P}\left(  \left[
n\right]  \right)  }\left(  -1\right)  ^{n-\left\vert I\right\vert }\left(
\sum_{i\in I}v_{i}\right)  ^{m}}+\sum_{I\in\mathcal{P}\left(  \left[
n+1\right]  \right)  \setminus\mathcal{P}\left(  \left[  n\right]  \right)
}\left(  -1\right)  ^{\left(  n+1\right)  -\left\vert I\right\vert }\left(
\sum_{i\in I}v_{i}\right)  ^{m}\\
&  =-\sum_{I\in\mathcal{P}\left(  \left[  n\right]  \right)  }\left(
-1\right)  ^{n-\left\vert I\right\vert }\left(  \sum_{i\in I}v_{i}\right)
^{m}+\sum_{I\in\mathcal{P}\left(  \left[  n+1\right]  \right)  \setminus
\mathcal{P}\left(  \left[  n\right]  \right)  }\left(  -1\right)  ^{\left(
n+1\right)  -\left\vert I\right\vert }\left(  \sum_{i\in I}v_{i}\right)  ^{m}.
\end{align*}
Adding $\sum_{I\in\mathcal{P}\left(  \left[  n\right]  \right)  }\left(
-1\right)  ^{n-\left\vert I\right\vert }\left(  \sum_{i\in I}v_{i}\right)
^{m}$ to both sides of this equality, we obtain%
\begin{align}
&  \sum_{I\in\mathcal{P}\left(  \left[  n\right]  \right)  }\left(  -1\right)
^{n-\left\vert I\right\vert }\left(  \sum_{i\in I}v_{i}\right)  ^{m}%
\nonumber\\
&  =\sum_{I\in\mathcal{P}\left(  \left[  n+1\right]  \right)  \setminus
\mathcal{P}\left(  \left[  n\right]  \right)  }\left(  -1\right)  ^{\left(
n+1\right)  -\left\vert I\right\vert }\left(  \sum_{i\in I}v_{i}\right)
^{m}\nonumber\\
&  =\sum_{U\in\mathcal{P}\left(  \left[  n\right]  \right)  }\left(
-1\right)  ^{\left(  n+1\right)  -\left\vert U\cup\left\{  n+1\right\}
\right\vert }\left(  \sum_{i\in U\cup\left\{  n+1\right\}  }v_{i}\right)
^{m}\nonumber\\
&  \ \ \ \ \ \ \ \ \ \ \left(
\begin{array}
[c]{c}%
\text{here, we have substituted }U\cup\left\{  n+1\right\}  \text{ for
}I\text{ in the sum,}\\
\text{since the map }\mathcal{P}\left(  \left[  n\right]  \right)
\rightarrow\mathcal{P}\left(  \left[  n+1\right]  \right)  \setminus
\mathcal{P}\left(  \left[  n\right]  \right)  ,\ U\mapsto U\cup\left\{
n+1\right\} \\
\text{is a bijection}%
\end{array}
\right) \nonumber\\
&  =\sum_{I\in\mathcal{P}\left(  \left[  n\right]  \right)  }\left(
-1\right)  ^{\left(  n+1\right)  -\left\vert I\cup\left\{  n+1\right\}
\right\vert }\left(  \sum_{i\in I\cup\left\{  n+1\right\}  }v_{i}\right)  ^{m}
\label{sol.noncomm.polarization2.wout.pf.1}%
\end{align}
(here, we have renamed the summation index $U$ as $I$).

But each $I\in\mathcal{P}\left(  \left[  n\right]  \right)  $ satisfies%
\begin{equation}
\left(  -1\right)  ^{\left(  n+1\right)  -\left\vert I\cup\left\{
n+1\right\}  \right\vert }\left(  \sum_{i\in I\cup\left\{  n+1\right\}  }%
v_{i}\right)  ^{m}=\left(  -1\right)  ^{n-\left\vert I\right\vert }\left(
w+\sum_{i\in I}v_{i}\right)  ^{m} \label{sol.noncomm.polarization2.wout.pf.2}%
\end{equation}
\footnote{\textit{Proof of (\ref{sol.noncomm.polarization2.wout.pf.2}):} Let
$I\in\mathcal{P}\left(  \left[  n\right]  \right)  $. But $\mathcal{P}\left(
\left[  n\right]  \right)  $ is the set of all subsets of $\left[  n\right]  $
(by the definition of $\mathcal{P}\left(  \left[  n\right]  \right)  $).
\par
Now, $I\in\mathcal{P}\left(  \left[  n\right]  \right)  $. Thus, $I$ is a
subset of $\left[  n\right]  $ (since $\mathcal{P}\left(  \left[  n\right]
\right)  $ is the set of all subsets of $\left[  n\right]  $). In other words,
$I\subseteq\left[  n\right]  $. If we had $n+1\in I$, then we would have
$n+1\in I\subseteq\left[  n\right]  =\left\{  1,2,\ldots,n\right\}  $, which
would contradict the fact that $n+1\notin\left\{  1,2,\ldots,n\right\}  $.
Thus, we cannot have $n+1\in I$. Hence, we have $n+1\notin I$. Hence,
$\left\vert I\cup\left\{  n+1\right\}  \right\vert =\left\vert I\right\vert
+1$ and $I\setminus\left\{  n+1\right\}  =I$. Therefore,
\[
\left(  n+1\right)  -\underbrace{\left\vert I\cup\left\{  n+1\right\}
\right\vert }_{=\left\vert I\right\vert +1}=\left(  n+1\right)  -\left(
\left\vert I\right\vert +1\right)  =n-\left\vert I\right\vert .
\]
Furthermore, we have $n+1\in\left\{  n+1\right\}  \subseteq I\cup\left\{
n+1\right\}  $. Hence, we can split off the addend for $i=n+1$ from the sum
$\sum_{i\in I\cup\left\{  n+1\right\}  }v_{i}$. We thus obtain%
\[
\sum_{i\in I\cup\left\{  n+1\right\}  }v_{i}=\underbrace{v_{n+1}}%
_{=w}+\underbrace{\sum_{\substack{i\in I\cup\left\{  n+1\right\}  ;\\i\neq
n+1}}}_{\substack{=\sum_{i\in\left(  I\cup\left\{  n+1\right\}  \right)
\setminus\left\{  n+1\right\}  }\\=\sum_{i\in I}\\\text{(since }\left(
I\cup\left\{  n+1\right\}  \right)  \setminus\left\{  n+1\right\}
=I\setminus\left\{  n+1\right\}  =I\text{)}}}v_{i}=w+\sum_{i\in I}v_{i}.
\]
Thus,%
\[
\underbrace{\left(  -1\right)  ^{\left(  n+1\right)  -\left\vert I\cup\left\{
n+1\right\}  \right\vert }}_{\substack{=\left(  -1\right)  ^{n-\left\vert
I\right\vert }\\\text{(since }\left(  n+1\right)  -\left\vert I\cup\left\{
n+1\right\}  \right\vert =n-\left\vert I\right\vert \text{)}}}\left(
\underbrace{\sum_{i\in I\cup\left\{  n+1\right\}  }v_{i}}_{=w+\sum_{i\in
I}v_{i}}\right)  ^{m}=\left(  -1\right)  ^{n-\left\vert I\right\vert }\left(
w+\sum_{i\in I}v_{i}\right)  ^{m}.
\]
This proves (\ref{sol.noncomm.polarization2.wout.pf.2}).}. Thus,
(\ref{sol.noncomm.polarization2.wout.pf.1}) becomes%
\begin{align*}
\sum_{I\in\mathcal{P}\left(  \left[  n\right]  \right)  }\left(  -1\right)
^{n-\left\vert I\right\vert }\left(  \sum_{i\in I}v_{i}\right)  ^{m}  &
=\sum_{I\in\mathcal{P}\left(  \left[  n\right]  \right)  }\underbrace{\left(
-1\right)  ^{\left(  n+1\right)  -\left\vert I\cup\left\{  n+1\right\}
\right\vert }\left(  \sum_{i\in I\cup\left\{  n+1\right\}  }v_{i}\right)
^{m}}_{\substack{=\left(  -1\right)  ^{n-\left\vert I\right\vert }\left(
w+\sum_{i\in I}v_{i}\right)  ^{m}\\\text{(by
(\ref{sol.noncomm.polarization2.wout.pf.2}))}}}\\
&  =\sum_{I\in\mathcal{P}\left(  \left[  n\right]  \right)  }\left(
-1\right)  ^{n-\left\vert I\right\vert }\left(  w+\sum_{i\in I}v_{i}\right)
^{m}.
\end{align*}
In view of $\sum_{I\in\mathcal{P}\left(  \left[  n\right]  \right)  }%
=\sum_{I\subseteq\left[  n\right]  }$ (an equality between summation signs),
this equality rewrites as%
\[
\sum_{I\subseteq\left[  n\right]  }\left(  -1\right)  ^{n-\left\vert
I\right\vert }\left(  \sum_{i\in I}v_{i}\right)  ^{m}=\sum_{I\subseteq\left[
n\right]  }\left(  -1\right)  ^{n-\left\vert I\right\vert }\left(
w+\sum_{i\in I}v_{i}\right)  ^{m}.
\]
In other words,%
\[
\sum_{I\subseteq\left[  n\right]  }\left(  -1\right)  ^{n-\left\vert
I\right\vert }\left(  w+\sum_{i\in I}v_{i}\right)  ^{m}=\sum_{I\subseteq
\left[  n\right]  }\left(  -1\right)  ^{n-\left\vert I\right\vert }\left(
\sum_{i\in I}v_{i}\right)  ^{m}.
\]
This proves (\ref{sol.noncomm.polarization2.wout}).]
\end{verlong}

\textbf{(a)} Let $m\in\left\{  0,1,\ldots,n-1\right\}  $ and $w\in\mathbb{L}$.
We have $m\in\left\{  0,1,\ldots,n-1\right\}  \subseteq\left\{  0,1,\ldots
,n\right\}  $. Hence, (\ref{sol.noncomm.polarization2.wout}) yields%
\[
\sum_{I\subseteq\left[  n\right]  }\left(  -1\right)  ^{n-\left\vert
I\right\vert }\left(  w+\sum_{i\in I}v_{i}\right)  ^{m}=\sum_{I\subseteq
\left[  n\right]  }\left(  -1\right)  ^{n-\left\vert I\right\vert }\left(
\sum_{i\in I}v_{i}\right)  ^{m}=0
\]
(by Exercise \ref{exe.noncomm.polarization} \textbf{(b)}). This solves
Exercise \ref{exe.noncomm.polarization2} \textbf{(a)}.

\textbf{(b)} Let $w\in\mathbb{L}$. We have $n\in\left\{  0,1,\ldots,n\right\}
$. Hence, (\ref{sol.noncomm.polarization2.wout}) (applied to $m=n$) yields%
\[
\sum_{I\subseteq\left[  n\right]  }\left(  -1\right)  ^{n-\left\vert
I\right\vert }\left(  w+\sum_{i\in I}v_{i}\right)  ^{n}=\sum_{I\subseteq
\left[  n\right]  }\left(  -1\right)  ^{n-\left\vert I\right\vert }\left(
\sum_{i\in I}v_{i}\right)  ^{n}=\sum_{\sigma\in S_{n}}v_{\sigma\left(
1\right)  }v_{\sigma\left(  2\right)  }\cdots v_{\sigma\left(  n\right)  }%
\]
(by Exercise \ref{exe.noncomm.polarization} \textbf{(c)}). This solves
Exercise \ref{exe.noncomm.polarization2} \textbf{(b)}.

Next, let us prepare for the solutions of parts \textbf{(c)} and \textbf{(d)}
of Exercise \ref{exe.noncomm.polarization2}. Set $q=\sum_{i\in\left[
n\right]  }v_{i}$. Then, each subset $I$ of $\left[  n\right]  $ satisfies%
\begin{align*}
q  &  =\sum_{i\in\left[  n\right]  }v_{i}=\underbrace{\sum_{\substack{i\in
\left[  n\right]  ;\\i\in I}}}_{\substack{=\sum_{i\in I}\\\text{(since
}I\subseteq\left[  n\right]  \text{)}}}v_{i}+\underbrace{\sum_{\substack{i\in
\left[  n\right]  ;\\i\notin I}}}_{=\sum_{i\in\left[  n\right]  \setminus I}%
}v_{i}\\
&  \ \ \ \ \ \ \ \ \ \ \left(  \text{since each }i\in\left[  n\right]  \text{
satisfies either }i\in I\text{ or }i\notin I\text{ (but not both)}\right) \\
&  =\sum_{i\in I}v_{i}+\sum_{i\in\left[  n\right]  \setminus I}v_{i}%
\end{align*}
and thus%
\begin{align}
-\underbrace{q}_{=\sum_{i\in I}v_{i}+\sum_{i\in\left[  n\right]  \setminus
I}v_{i}}+\underbrace{\sum_{i\in I}2v_{i}}_{=2\sum_{i\in I}v_{i}}  &  =-\left(
\sum_{i\in I}v_{i}+\sum_{i\in\left[  n\right]  \setminus I}v_{i}\right)
+2\sum_{i\in I}v_{i}\nonumber\\
&  =\sum_{i\in I}v_{i}-\sum_{i\in\left[  n\right]  \setminus I}v_{i}.
\label{sol.noncomm.polarization2.-q}%
\end{align}

\textbf{(c)} Let $m\in\left\{  0,1,\ldots,n-1\right\}  $. Then, Exercise
\ref{exe.noncomm.polarization2} \textbf{(a)} (applied to $2v_{i}$ and $-q$
instead of $v_{i}$ and $w$) yields%
\[
\sum_{I\subseteq\left[  n\right]  }\left(  -1\right)  ^{n-\left\vert
I\right\vert }\left(  -q+\sum_{i\in I}2v_{i}\right)  ^{m}=0.
\]

\begin{vershort}
\noindent Using (\ref{sol.noncomm.polarization2.-q}), we can rewrite this as%
\[
\sum_{I\subseteq\left[  n\right]  }\left(  -1\right)  ^{n-\left\vert
I\right\vert }\left(  \sum_{i\in I}v_{i}-\sum_{i\in\left[  n\right]  \setminus
I}v_{i}\right)  ^{m}=0.
\]

\end{vershort}

\begin{verlong}
\noindent Comparing this with%
\[
\sum_{I\subseteq\left[  n\right]  }\left(  -1\right)  ^{n-\left\vert
I\right\vert }\left(  \underbrace{-q+\sum_{i\in I}2v_{i}}_{\substack{=\sum
_{i\in I}v_{i}-\sum_{i\in\left[  n\right]  \setminus I}v_{i}\\\text{(by
(\ref{sol.noncomm.polarization2.-q}))}}}\right)  ^{m}=\sum_{I\subseteq\left[
n\right]  }\left(  -1\right)  ^{n-\left\vert I\right\vert }\left(  \sum_{i\in
I}v_{i}-\sum_{i\in\left[  n\right]  \setminus I}v_{i}\right)  ^{m},
\]
we obtain%
\[
\sum_{I\subseteq\left[  n\right]  }\left(  -1\right)  ^{n-\left\vert
I\right\vert }\left(  \sum_{i\in I}v_{i}-\sum_{i\in\left[  n\right]  \setminus
I}v_{i}\right)  ^{m}=0.
\]

\end{verlong}

\noindent This solves Exercise \ref{exe.noncomm.polarization2} \textbf{(c)}.

\begin{vershort}
\textbf{(d)} Exercise \ref{exe.noncomm.polarization2} \textbf{(b)} (applied to
$2v_{i}$ and $-q$ instead of $v_{i}$ and $w$) yields%
\begin{align*}
\sum_{I\subseteq\left[  n\right]  }\left(  -1\right)  ^{n-\left\vert
I\right\vert }\left(  -q+\sum_{i\in I}2v_{i}\right)  ^{n}  &  =\sum_{\sigma\in
S_{n}}\underbrace{\left(  2v_{\sigma\left(  1\right)  }\right)  \left(
2v_{\sigma\left(  2\right)  }\right)  \cdots\left(  2v_{\sigma\left(
n\right)  }\right)  }_{=2^{n}\cdot v_{\sigma\left(  1\right)  }v_{\sigma
\left(  2\right)  }\cdots v_{\sigma\left(  n\right)  }}\\
&  =2^{n}\sum_{\sigma\in S_{n}}v_{\sigma\left(  1\right)  }v_{\sigma\left(
2\right)  }\cdots v_{\sigma\left(  n\right)  }.
\end{align*}
Using (\ref{sol.noncomm.polarization2.-q}), we can rewrite this as%
\[
\sum_{I\subseteq\left[  n\right]  }\left(  -1\right)  ^{n-\left\vert
I\right\vert }\left(  \sum_{i\in I}v_{i}-\sum_{i\in\left[  n\right]  \setminus
I}v_{i}\right)  ^{n}=2^{n}\sum_{\sigma\in S_{n}}v_{\sigma\left(  1\right)
}v_{\sigma\left(  2\right)  }\cdots v_{\sigma\left(  n\right)  }.
\]
This solves Exercise \ref{exe.noncomm.polarization2} \textbf{(d)}. \qedhere

\end{vershort}

\begin{verlong}
\textbf{(d)} Exercise \ref{exe.noncomm.polarization2} \textbf{(b)} (applied to
$2v_{i}$ and $-q$ instead of $v_{i}$ and $w$) yields%
\begin{align*}
\sum_{I\subseteq\left[  n\right]  }\left(  -1\right)  ^{n-\left\vert
I\right\vert }\left(  -q+\sum_{i\in I}2v_{i}\right)  ^{n}  &  =\sum_{\sigma\in
S_{n}}\underbrace{\left(  2v_{\sigma\left(  1\right)  }\right)  \left(
2v_{\sigma\left(  2\right)  }\right)  \cdots\left(  2v_{\sigma\left(
n\right)  }\right)  }_{=2^{n}\cdot v_{\sigma\left(  1\right)  }v_{\sigma
\left(  2\right)  }\cdots v_{\sigma\left(  n\right)  }}\\
&  =\sum_{\sigma\in S_{n}}2^{n}\cdot v_{\sigma\left(  1\right)  }%
v_{\sigma\left(  2\right)  }\cdots v_{\sigma\left(  n\right)  }\\
&  =2^{n}\sum_{\sigma\in S_{n}}v_{\sigma\left(  1\right)  }v_{\sigma\left(
2\right)  }\cdots v_{\sigma\left(  n\right)  }.
\end{align*}
Comparing this with%
\[
\sum_{I\subseteq\left[  n\right]  }\left(  -1\right)  ^{n-\left\vert
I\right\vert }\left(  \underbrace{-q+\sum_{i\in I}2v_{i}}_{\substack{=\sum
_{i\in I}v_{i}-\sum_{i\in\left[  n\right]  \setminus I}v_{i}\\\text{(by
(\ref{sol.noncomm.polarization2.-q}))}}}\right)  ^{n}=\sum_{I\subseteq\left[
n\right]  }\left(  -1\right)  ^{n-\left\vert I\right\vert }\left(  \sum_{i\in
I}v_{i}-\sum_{i\in\left[  n\right]  \setminus I}v_{i}\right)  ^{n},
\]
we obtain%
\[
\sum_{I\subseteq\left[  n\right]  }\left(  -1\right)  ^{n-\left\vert
I\right\vert }\left(  \sum_{i\in I}v_{i}-\sum_{i\in\left[  n\right]  \setminus
I}v_{i}\right)  ^{n}=2^{n}\sum_{\sigma\in S_{n}}v_{\sigma\left(  1\right)
}v_{\sigma\left(  2\right)  }\cdots v_{\sigma\left(  n\right)  }.
\]
This solves Exercise \ref{exe.noncomm.polarization2} \textbf{(d)}.
\end{verlong}
\end{proof}

\subsection{Solution to Exercise \ref{exe.det.sumdets1}}

In this section, we shall use the same notations that we have introduced in
Section \ref{sect.sol.noncomm.polarization}.

We start by proving a simple restatement of Lemma \ref{lem.noncomm.prodrule2}:

\begin{lemma}
\label{lem.noncomm.prodrule2b}Let $\mathbb{L}$ be a noncommutative ring.

Let $n\in\mathbb{N}$. Let $I$ be a finite set. For every $i\in I$ and every
$j\in\left[  n\right]  $, we let $b_{i,j}$ be an element of $\mathbb{L}$.
Then,%
\[
\prod_{j=1}^{n}\sum_{i\in I}b_{i,j}=\sum_{f:\left[  n\right]  \rightarrow
I}\prod_{j=1}^{n}b_{f\left(  j\right)  ,j}.
\]

\end{lemma}

\begin{proof}
[Proof of Lemma \ref{lem.noncomm.prodrule2b}.]Define an $m\in\mathbb{N}$ by
$m=\left\vert I\right\vert $. (This is well-defined, since $I$ is a finite
set.) Now, $\left\vert \left[  m\right]  \right\vert =m=\left\vert
I\right\vert $. Hence, there exists a bijection $\varphi:\left[  m\right]
\rightarrow I$. Consider such a $\varphi$.

Lemma \ref{lem.noncomm.prodrule2} (applied to $\mathbb{K}=\mathbb{L}$ and
$p_{i,k}=b_{\varphi\left(  k\right)  ,i}$) yields%
\begin{equation}
\prod_{i=1}^{n}\sum_{k=1}^{m}b_{\varphi\left(  k\right)  ,i}=\sum
_{\kappa:\left[  n\right]  \rightarrow\left[  m\right]  }\prod_{i=1}%
^{n}b_{\varphi\left(  \kappa\left(  i\right)  \right)  ,i}.
\label{pf.lem.noncomm.prodrule2b.0}%
\end{equation}

Every $i\in\left\{  1,2,\ldots,n\right\}  $ satisfies
\begin{equation}
\sum_{k=1}^{m}b_{\varphi\left(  k\right)  ,i}=\sum_{t\in I}b_{t,i}
\label{pf.lem.noncomm.prodrule2b.1}%
\end{equation}
\footnote{\textit{Proof of (\ref{pf.lem.noncomm.prodrule2b.1}):} Let
$i\in\left\{  1,2,\ldots,n\right\}  $. Recall that $\varphi:\left[  m\right]
\rightarrow I$ is a bijection. Hence, we can substitute $\varphi\left(
k\right)  $ for $t$ in the sum $\sum_{t\in I}b_{t,i}$. We thus obtain%
\[
\sum_{t\in I}b_{t,i}=\underbrace{\sum_{k\in\left[  m\right]  }}%
_{\substack{=\sum_{k\in\left\{  1,2,\ldots,m\right\}  }\\\text{(since }\left[
m\right]  =\left\{  1,2,\ldots,m\right\}  \text{)}}}b_{\varphi\left(
k\right)  ,i}=\underbrace{\sum_{k\in\left\{  1,2,\ldots,m\right\}  }}%
_{=\sum_{k=1}^{m}}b_{\varphi\left(  k\right)  ,i}=\sum_{k=1}^{m}%
b_{\varphi\left(  k\right)  ,i}.
\]
This proves (\ref{pf.lem.noncomm.prodrule2b.1}).}.

Hence,%
\begin{align}
\prod_{i=1}^{n}\underbrace{\sum_{k=1}^{m}b_{\varphi\left(  k\right)  ,i}%
}_{\substack{=\sum_{t\in I}b_{t,i}\\\text{(by
(\ref{pf.lem.noncomm.prodrule2b.1}))}}}  &  =\prod_{i=1}^{n}\sum_{t\in
I}b_{t,i}=\prod_{j=1}^{n}\underbrace{\sum_{t\in I}b_{t,j}}_{\substack{=\sum
_{i\in I}b_{i,j}\\\text{(here, we have renamed the}\\\text{summation index
}t\text{ as }i\text{)}}}\nonumber\\
&  \ \ \ \ \ \ \ \ \ \ \left(
\begin{array}
[c]{c}%
\text{here, we have renamed the index }i\text{ as }j\\
\text{in the product}%
\end{array}
\right) \nonumber\\
&  =\prod_{j=1}^{n}\sum_{i\in I}b_{i,j}. \label{pf.lem.noncomm.prodrule2b.L}%
\end{align}

Comparing this with (\ref{pf.lem.noncomm.prodrule2b.0}), we find%
\begin{equation}
\prod_{j=1}^{n}\sum_{i\in I}b_{i,j}=\sum_{\kappa:\left[  n\right]
\rightarrow\left[  m\right]  }\prod_{i=1}^{n}b_{\varphi\left(  \kappa\left(
i\right)  \right)  ,i}. \label{pf.lem.noncomm.prodrule2b.4}%
\end{equation}

\begin{vershort}
The map $\varphi$ is a bijection, and thus its inverse $\varphi^{-1}$ exists.
Hence, the map%
\[
\left[  m\right]  ^{\left[  n\right]  }\rightarrow I^{\left[  n\right]
},\ \ \ \ \ \ \ \ \ \ \kappa\mapsto\varphi\circ\kappa
\]
is a bijection (in fact, its inverse is the map $I^{\left[  n\right]
}\rightarrow\left[  m\right]  ^{\left[  n\right]  },\ f\mapsto\varphi
^{-1}\circ f$).
\end{vershort}

\begin{verlong}
The map $\varphi$ is a bijection, and thus is invertible. Hence, its inverse
$\varphi^{-1}:I\rightarrow\left[  m\right]  $ exists. Let $A$ be the map%
\[
\left[  m\right]  ^{\left[  n\right]  }\rightarrow I^{\left[  n\right]
},\ \ \ \ \ \ \ \ \ \ \kappa\mapsto\varphi\circ\kappa.
\]
Let $B$ be the map%
\[
I^{\left[  n\right]  }\rightarrow\left[  m\right]  ^{\left[  n\right]
},\ \ \ \ \ \ \ \ \ \ f\mapsto\varphi^{-1}\circ f.
\]
Then, $A\circ B=\operatorname*{id}$\ \ \ \ \footnote{\textit{Proof.} Let $f\in
I^{\left[  n\right]  }$ be arbitrary. The definition of $B$ yields $B\left(
f\right)  =\varphi^{-1}\circ f$. The definition of $A$ yields%
\[
A\left(  B\left(  f\right)  \right)  =\varphi\circ\underbrace{\left(  B\left(
f\right)  \right)  }_{=\varphi^{-1}\circ f}=\underbrace{\varphi\circ
\varphi^{-1}}_{=\operatorname*{id}}\circ f=f=\operatorname*{id}\left(
f\right)  .
\]
Thus, $\left(  A\circ B\right)  \left(  f\right)  =A\left(  B\left(  f\right)
\right)  =\operatorname*{id}\left(  f\right)  $.
\par
Now, forget that we fixed $f$. We thus have shown that $\left(  A\circ
B\right)  \left(  f\right)  =\operatorname*{id}\left(  f\right)  $ for each
$f\in I^{\left[  n\right]  }$. In other words, $A\circ B=\operatorname*{id}$.}
and $B\circ A=\operatorname*{id}$\ \ \ \ \footnote{\textit{Proof.} Let
$\kappa\in\left[  m\right]  ^{\left[  n\right]  }$ be arbitrary. The
definition of $A$ yields $A\left(  \kappa\right)  =\varphi\circ\kappa$. The
definition of $B$ yields%
\[
B\left(  A\left(  \kappa\right)  \right)  =\varphi^{-1}\circ
\underbrace{\left(  A\left(  \kappa\right)  \right)  }_{=\varphi\circ\kappa
}=\underbrace{\varphi^{-1}\circ\varphi}_{=\operatorname*{id}}\circ
\kappa=\kappa=\operatorname*{id}\left(  \kappa\right)  .
\]
Thus, $\left(  B\circ A\right)  \left(  \kappa\right)  =B\left(  A\left(
\kappa\right)  \right)  =\operatorname*{id}\left(  \kappa\right)  $.
\par
Now, forget that we fixed $\kappa$. We thus have shown that $\left(  B\circ
A\right)  \left(  \kappa\right)  =\operatorname*{id}\left(  \kappa\right)  $
for each $\kappa\in\left[  m\right]  ^{\left[  n\right]  }$. In other words,
$B\circ A=\operatorname*{id}$.}. Hence, the maps $A$ and $B$ are mutually
inverse. Thus, the map $A$ is invertible, and thus a bijection. In other
words, the map
\[
\left[  m\right]  ^{\left[  n\right]  }\rightarrow I^{\left[  n\right]
},\ \ \ \ \ \ \ \ \ \ \kappa\mapsto\varphi\circ\kappa
\]
is a bijection\footnote{since the map%
\[
\left[  m\right]  ^{\left[  n\right]  }\rightarrow I^{\left[  n\right]
},\ \ \ \ \ \ \ \ \ \ \kappa\mapsto\varphi\circ\kappa
\]
is precisely the map $A$}.
\end{verlong}

Therefore, we can substitute $\varphi\circ\kappa$ for $f$ in the sum
$\sum_{f\in I^{\left[  n\right]  }}\prod_{j=1}^{n}b_{f\left(  j\right)  ,j}$.
We thus obtain%
\begin{align*}
\sum_{f\in I^{\left[  n\right]  }}\prod_{j=1}^{n}b_{f\left(  j\right)  ,j}  &
=\underbrace{\sum_{\kappa\in\left[  m\right]  ^{\left[  n\right]  }}}%
_{=\sum_{\kappa:\left[  n\right]  \rightarrow\left[  m\right]  }}\prod
_{j=1}^{n}\underbrace{b_{\left(  \varphi\circ\kappa\right)  \left(  j\right)
,j}}_{\substack{=b_{\varphi\left(  \kappa\left(  j\right)  \right)
,j}\\\text{(since }\left(  \varphi\circ\kappa\right)  \left(  j\right)
=\varphi\left(  \kappa\left(  j\right)  \right)  \text{)}}}=\sum
_{\kappa:\left[  n\right]  \rightarrow\left[  m\right]  }\underbrace{\prod
_{j=1}^{n}b_{\varphi\left(  \kappa\left(  j\right)  \right)  ,j}%
}_{\substack{=\prod_{i=1}^{n}b_{\varphi\left(  \kappa\left(  i\right)
\right)  ,i}\\\text{(here, we have renamed}\\\text{the index }j\text{ as
}i\text{ in}\\\text{the product)}}}\\
&  =\sum_{\kappa:\left[  n\right]  \rightarrow\left[  m\right]  }\prod
_{i=1}^{n}b_{\varphi\left(  \kappa\left(  i\right)  \right)  ,i}.
\end{align*}
Comparing this with (\ref{pf.lem.noncomm.prodrule2b.4}), we obtain%
\[
\prod_{j=1}^{n}\sum_{i\in I}b_{i,j}=\underbrace{\sum_{f\in I^{\left[
n\right]  }}}_{=\sum_{f:\left[  n\right]  \rightarrow I}}\prod_{j=1}%
^{n}b_{f\left(  j\right)  ,j}=\sum_{f:\left[  n\right]  \rightarrow I}%
\prod_{j=1}^{n}b_{f\left(  j\right)  ,j}.
\]
This proves Lemma \ref{lem.noncomm.prodrule2b}.
\end{proof}

Next, we state an analogue of Lemma \ref{lem.sol.noncomm.polarization.1}:

\begin{lemma}
\label{lem.sol.det.sumdets1.1}Let $\mathbb{L}$ be a noncommutative ring. Let
$G$ be a finite set. Let $n\in\mathbb{N}$.

For each $i\in G$ and $j\in\left[  n\right]  $, let $b_{i,j}$ be an element of
$\mathbb{L}$. Let $I$ be a subset of $G$. Then,%
\[
\prod_{j=1}^{n}\sum_{i\in I}b_{i,j}=\sum_{\substack{f:\left[  n\right]
\rightarrow G;\\f\left(  \left[  n\right]  \right)  \subseteq I}}\prod
_{j=1}^{n}b_{f\left(  j\right)  ,j}.
\]

\end{lemma}

\begin{proof}
[Proof of Lemma \ref{lem.sol.det.sumdets1.1}.]For each $j\in\left[  n\right]
$, we have%
\begin{align}
&  \sum_{i\in G}\left[  i\in I\right]  b_{i,j}\nonumber\\
&  =\underbrace{\sum_{\substack{i\in G;\\i\in I}}}_{\substack{=\sum_{i\in
I}\\\text{(since }I\subseteq G\text{)}}}\underbrace{\left[  i\in I\right]
}_{\substack{=1\\\text{(since }i\in I\text{ is true)}}}b_{i,j}+\sum
_{\substack{i\in G;\\i\notin I}}\underbrace{\left[  i\in I\right]
}_{\substack{=0\\\text{(since }i\in I\text{ is false}\\\text{(because }i\notin
I\text{))}}}b_{i,j}\nonumber\\
&  \ \ \ \ \ \ \ \ \ \ \left(  \text{since each }i\in G\text{ satisfies either
}i\in I\text{ or }i\notin I\text{ (but not both)}\right) \nonumber\\
&  =\sum_{i\in I}1b_{i,j}+\underbrace{\sum_{\substack{i\in G;\\i\notin
I}}0b_{i,j}}_{=0}=\sum_{i\in I}1b_{i,j}=\sum_{i\in I}b_{i,j}.
\label{pf.lem.sol.det.sumdets1.1.1}%
\end{align}
Hence,%
\[
\prod_{j=1}^{n}\underbrace{\sum_{i\in G}\left[  i\in I\right]  b_{i,j}%
}_{\substack{=\sum_{i\in I}b_{i,j}\\\text{(by
(\ref{pf.lem.sol.det.sumdets1.1.1}))}}}=\prod_{j=1}^{n}\sum_{i\in I}b_{i,j}.
\]
Thus,%
\begin{equation}
\prod_{j=1}^{n}\sum_{i\in I}b_{i,j}=\prod_{j=1}^{n}\sum_{i\in G}\left[  i\in
I\right]  b_{i,j}=\sum_{f:\left[  n\right]  \rightarrow G}\prod_{j=1}%
^{n}\left(  \left[  f\left(  j\right)  \in I\right]  b_{f\left(  j\right)
,j}\right)  \label{pf.lem.sol.det.sumdets1.1.4}%
\end{equation}
(by Lemma \ref{lem.noncomm.prodrule2b} (applied to $G$ and $\left[  i\in
I\right]  b_{i,j}$ instead of $I$ and $b_{i,j}$)).

Now, fix any map $f:\left[  n\right]  \rightarrow G$. Then,%
\begin{align}
&  \prod_{j=1}^{n}\left(  \left[  f\left(  j\right)  \in I\right]  b_{f\left(
j\right)  ,j}\right) \nonumber\\
&  =\left(  \left[  f\left(  1\right)  \in I\right]  b_{f\left(  1\right)
,1}\right)  \left(  \left[  f\left(  2\right)  \in I\right]  b_{f\left(
2\right)  ,2}\right)  \cdots\left(  \left[  f\left(  n\right)  \in I\right]
b_{f\left(  n\right)  ,n}\right)  . \label{pf.lem.sol.det.sumdets1.1.6}%
\end{align}
The factors $\left[  f\left(  1\right)  \in I\right]  ,\left[  f\left(
2\right)  \in I\right]  ,\ldots,\left[  f\left(  n\right)  \in I\right]  $ on
the right hand side of this equality are integers, and therefore can be freely
moved within the product (even though $\mathbb{L}$ is not necessarily
commutative). In particular, we can move them to the front of the product.
Thus, we find%
\begin{align*}
&  \left(  \left[  f\left(  1\right)  \in I\right]  b_{f\left(  1\right)
,1}\right)  \left(  \left[  f\left(  2\right)  \in I\right]  b_{f\left(
2\right)  ,2}\right)  \cdots\left(  \left[  f\left(  n\right)  \in I\right]
b_{f\left(  n\right)  ,n}\right) \\
&  =\underbrace{\left(  \left[  f\left(  1\right)  \in I\right]  \left[
f\left(  2\right)  \in I\right]  \cdots\left[  f\left(  n\right)  \in
I\right]  \right)  }_{=\prod_{i=1}^{n}\left[  f\left(  i\right)  \in I\right]
}\left(  b_{f\left(  1\right)  ,1}b_{f\left(  2\right)  ,2}\cdots b_{f\left(
n\right)  ,n}\right) \\
&  =\left(  \prod_{i=1}^{n}\left[  f\left(  i\right)  \in I\right]  \right)
\left(  b_{f\left(  1\right)  ,1}b_{f\left(  2\right)  ,2}\cdots b_{f\left(
n\right)  ,n}\right)  .
\end{align*}
Hence, (\ref{pf.lem.sol.det.sumdets1.1.6}) rewrites as%
\begin{equation}
\prod_{j=1}^{n}\left(  \left[  f\left(  j\right)  \in I\right]  b_{f\left(
j\right)  ,j}\right)  =\left(  \prod_{i=1}^{n}\left[  f\left(  i\right)  \in
I\right]  \right)  \left(  b_{f\left(  1\right)  ,1}b_{f\left(  2\right)
,2}\cdots b_{f\left(  n\right)  ,n}\right)  .
\label{pf.lem.sol.det.sumdets1.1.6b}%
\end{equation}

But Lemma \ref{lem.sol.noncomm.polarization.1.7} (applied to $n$ instead of
$m$) yields%
\begin{equation}
\prod_{i=1}^{n}\left[  f\left(  i\right)  \in I\right]  =\left[  f\left(
\left[  n\right]  \right)  \subseteq I\right]  .
\label{pf.lem.sol.det.sumdets1.1.7n}%
\end{equation}

Hence, (\ref{pf.lem.sol.det.sumdets1.1.6b}) becomes%
\begin{align}
\prod_{j=1}^{n}\left(  \left[  f\left(  j\right)  \in I\right]  b_{f\left(
j\right)  ,j}\right)   &  =\underbrace{\left(  \prod_{i=1}^{n}\left[  f\left(
i\right)  \in I\right]  \right)  }_{\substack{=\left[  f\left(  \left[
n\right]  \right)  \subseteq I\right]  \\\text{(by
(\ref{pf.lem.sol.det.sumdets1.1.7n}))}}}\left(  b_{f\left(  1\right)
,1}b_{f\left(  2\right)  ,2}\cdots b_{f\left(  n\right)  ,n}\right)
\nonumber\\
&  =\left[  f\left(  \left[  n\right]  \right)  \subseteq I\right]  \left(
b_{f\left(  1\right)  ,1}b_{f\left(  2\right)  ,2}\cdots b_{f\left(  n\right)
,n}\right)  . \label{pf.lem.sol.det.sumdets1.1.8}%
\end{align}

Now, forget that we fixed $f$. We thus have proven
(\ref{pf.lem.sol.det.sumdets1.1.8}) for each map $f:\left[  n\right]
\rightarrow G$.

Now, (\ref{pf.lem.sol.det.sumdets1.1.4}) becomes%
\begin{align*}
\prod_{j=1}^{n}\sum_{i\in I}b_{i,j}  &  =\sum_{f:\left[  n\right]  \rightarrow
G}\underbrace{\prod_{j=1}^{n}\left(  \left[  f\left(  j\right)  \in I\right]
b_{f\left(  j\right)  ,j}\right)  }_{\substack{=\left[  f\left(  \left[
n\right]  \right)  \subseteq I\right]  \left(  b_{f\left(  1\right)
,1}b_{f\left(  2\right)  ,2}\cdots b_{f\left(  n\right)  ,n}\right)
\\\text{(by (\ref{pf.lem.sol.det.sumdets1.1.8}))}}}\\
&  =\sum_{f:\left[  n\right]  \rightarrow G}\left[  f\left(  \left[  n\right]
\right)  \subseteq I\right]  \left(  b_{f\left(  1\right)  ,1}b_{f\left(
2\right)  ,2}\cdots b_{f\left(  n\right)  ,n}\right) \\
&  =\sum_{\substack{f:\left[  n\right]  \rightarrow G;\\f\left(  \left[
n\right]  \right)  \subseteq I}}\underbrace{\left[  f\left(  \left[  n\right]
\right)  \subseteq I\right]  }_{\substack{=1\\\text{(since }f\left(  \left[
n\right]  \right)  \subseteq I\text{)}}}\left(  b_{f\left(  1\right)
,1}b_{f\left(  2\right)  ,2}\cdots b_{f\left(  n\right)  ,n}\right) \\
&  \ \ \ \ \ \ \ \ \ \ +\sum_{\substack{f:\left[  n\right]  \rightarrow
G;\\\text{not }f\left(  \left[  n\right]  \right)  \subseteq I}%
}\underbrace{\left[  f\left(  \left[  n\right]  \right)  \subseteq I\right]
}_{\substack{=0\\\text{(since we don't have }f\left(  \left[  n\right]
\right)  \subseteq I\text{)}}}\left(  b_{f\left(  1\right)  ,1}b_{f\left(
2\right)  ,2}\cdots b_{f\left(  n\right)  ,n}\right) \\
&  \ \ \ \ \ \ \ \ \ \ \left(
\begin{array}
[c]{c}%
\text{since each map }f:\left[  n\right]  \rightarrow G\text{ satisfies either
}f\left(  \left[  n\right]  \right)  \subseteq I\\
\text{or }\left(  \text{not }f\left(  \left[  n\right]  \right)  \subseteq
I\right)  \text{ (but not both)}%
\end{array}
\right) \\
&  =\sum_{\substack{f:\left[  n\right]  \rightarrow G;\\f\left(  \left[
n\right]  \right)  \subseteq I}}b_{f\left(  1\right)  ,1}b_{f\left(  2\right)
,2}\cdots b_{f\left(  n\right)  ,n}+\underbrace{\sum_{\substack{f:\left[
n\right]  \rightarrow G;\\\text{not }f\left(  \left[  n\right]  \right)
\subseteq I}}0\left(  b_{f\left(  1\right)  ,1}b_{f\left(  2\right)  ,2}\cdots
b_{f\left(  n\right)  ,n}\right)  }_{=0}\\
&  =\sum_{\substack{f:\left[  n\right]  \rightarrow G;\\f\left(  \left[
n\right]  \right)  \subseteq I}}\underbrace{b_{f\left(  1\right)
,1}b_{f\left(  2\right)  ,2}\cdots b_{f\left(  n\right)  ,n}}_{=\prod
_{j=1}^{n}b_{f\left(  j\right)  ,j}}=\sum_{\substack{f:\left[  n\right]
\rightarrow G;\\f\left(  \left[  n\right]  \right)  \subseteq I}}\prod
_{j=1}^{n}b_{f\left(  j\right)  ,j}.
\end{align*}
This proves Lemma \ref{lem.sol.det.sumdets1.1}.
\end{proof}

The next result we shall use can be viewed as a slight generalization of Lemma
\ref{lem.sol.noncomm.polarization.2}:

\begin{lemma}
\label{lem.sol.det.sumdets1.2}Let $G$ be a finite set. Let $H$ and $T$ be two
subsets of $G$. Then,%
\[
\sum_{\substack{I\subseteq G;\\H\subseteq I;\\T\subseteq I}}\left(  -1\right)
^{\left\vert I\right\vert }=\left(  -1\right)  ^{\left\vert G\right\vert
}\left[  G\setminus H\subseteq T\right]  .
\]

\end{lemma}

\begin{proof}
[Proof of Lemma \ref{lem.sol.det.sumdets1.2}.]If $I$ is a subset of $G$, then
the statement $\left(  H\subseteq I\text{ and }T\subseteq I\right)  $ is
equivalent to the statement $\left(  H\cup T\subseteq I\right)  $%
\ \ \ \ \footnote{In fact, this equivalence is a basic fact of set theory,
which holds for any set $I$.}. Hence, we have the following equality of
summation signs:%
\begin{equation}
\sum_{\substack{I\subseteq G;\\H\subseteq I;\\T\subseteq I}}=\sum
_{\substack{I\subseteq G;\\H\cup T\subseteq I}}.
\label{pf.lem.sol.det.sumdets1.2.sumeq}%
\end{equation}

Let $R=H\cup T$. Then, $R$ is a subset of $G$ (since $H$ and $T$ are subsets
of $G$). Moreover, the equality (\ref{pf.lem.sol.det.sumdets1.2.sumeq})
becomes%
\begin{equation}
\sum_{\substack{I\subseteq G;\\H\subseteq I;\\T\subseteq I}}=\sum
_{\substack{I\subseteq G;\\H\cup T\subseteq I}}=\sum_{\substack{I\subseteq
G;\\R\subseteq I}} \label{pf.lem.sol.det.sumdets1.2.sumeq2}%
\end{equation}
(since $H\cup T=R$).

We are in one of the following two cases:

\textit{Case 1:} We have $G\setminus H\not \subseteq T$.

\textit{Case 2:} We have $G\setminus H\subseteq T$.

Let us first consider Case 1. In this case, we have $G\setminus
H\not \subseteq T$. Hence, there exists some $g\in G\setminus H$ such that
$g\notin T$. Consider this $g$.

We have $g\in G\setminus H$. In other words, $g\in G$ and $g\notin H$. We have
$g\notin H\cup T$ (since $g\notin H$ and $g\notin T$). In other words,
$g\notin R$ (since $R=H\cup T$). Also, $\left\{  g\right\}  $ is a subset of
$G$ (since $g\in G$).

\begin{vershort}
Now, the map%
\begin{align}
\left\{  I\subseteq G\ \mid\ R\subseteq I\text{ and }g\in I\right\}   &
\rightarrow\left\{  I\subseteq G\ \mid\ R\subseteq I\text{ and }g\notin
I\right\}  ,\nonumber\\
J  &  \mapsto J\setminus\left\{  g\right\}
\label{pf.lem.sol.det.sumdets1.2.short.14}%
\end{align}
is well-defined\footnote{This is because every $J\in\left\{  I\subseteq
G\ \mid\ R\subseteq I\text{ and }g\in I\right\}  $ satisfies $J\setminus
\left\{  g\right\}  \in\left\{  I\subseteq G\ \mid\ R\subseteq I\text{ and
}g\notin I\right\}  $. (Proving this is straightforward using the facts that
$g\in G$ and $g\notin R$.)} and bijective\footnote{Indeed, its inverse is the
map%
\begin{align*}
\left\{  I\subseteq G\ \mid\ R\subseteq I\text{ and }g\notin I\right\}   &
\rightarrow\left\{  I\subseteq G\ \mid\ R\subseteq I\text{ and }g\in
I\right\}  ,\\
J  &  \mapsto J\cup\left\{  g\right\}
\end{align*}
(which, again, is well-defined because $g\in G$).}.

Hence,%
\begin{align}
&  \underbrace{\sum_{\substack{J\subseteq G;\\R\subseteq J;\\g\notin J}%
}}_{=\sum_{J\in\left\{  I\subseteq G\ \mid\ R\subseteq I\text{ and }g\notin
I\right\}  }}\left(  -1\right)  ^{\left\vert J\right\vert }\nonumber\\
&  =\sum_{J\in\left\{  I\subseteq G\ \mid\ R\subseteq I\text{ and }g\notin
I\right\}  }\left(  -1\right)  ^{\left\vert J\right\vert }=\underbrace{\sum
_{J\in\left\{  I\subseteq G\ \mid\ R\subseteq I\text{ and }g\in I\right\}  }%
}_{=\sum_{\substack{J\subseteq G;\\R\subseteq J;\\g\in J}}}\left(  -1\right)
^{\left\vert J\setminus\left\{  g\right\}  \right\vert }\nonumber\\
&  \ \ \ \ \ \ \ \ \ \ \left(
\begin{array}
[c]{c}%
\text{here, we have substituted }J\setminus\left\{  g\right\}  \text{ for
}J\text{ in the sum,}\\
\text{since the map (\ref{pf.lem.sol.det.sumdets1.2.short.14}) is bijective}%
\end{array}
\right) \nonumber\\
&  =\sum_{\substack{J\subseteq G;\\R\subseteq J;\\g\in J}}\underbrace{\left(
-1\right)  ^{\left\vert J\setminus\left\{  g\right\}  \right\vert }%
}_{\substack{=\left(  -1\right)  ^{\left\vert J\right\vert -1}\\\text{(since
}\left\vert J\setminus\left\{  g\right\}  \right\vert =\left\vert J\right\vert
-1\\\text{(because }g\in J\text{))}}}=\sum_{\substack{J\subseteq
G;\\R\subseteq J;\\g\in J}}\underbrace{\left(  -1\right)  ^{\left\vert
J\right\vert -1}}_{=-\left(  -1\right)  ^{\left\vert J\right\vert }}%
=\sum_{\substack{J\subseteq G;\\R\subseteq J;\\g\in J}}\left(  -\left(
-1\right)  ^{\left\vert J\right\vert }\right) \nonumber\\
&  =-\sum_{\substack{J\subseteq G;\\R\subseteq J;\\g\in J}}\left(  -1\right)
^{\left\vert J\right\vert }. \label{pf.lem.sol.det.sumdets1.2.short.15}%
\end{align}

But (\ref{pf.lem.sol.det.sumdets1.2.sumeq2}) yields
\begin{align*}
\sum_{\substack{I\subseteq G;\\H\subseteq I;\\T\subseteq I}}\left(  -1\right)
^{\left\vert I\right\vert }  &  =\sum_{\substack{I\subseteq G;\\R\subseteq
I}}\left(  -1\right)  ^{\left\vert I\right\vert }=\sum_{\substack{J\subseteq
G;\\R\subseteq J}}\left(  -1\right)  ^{\left\vert J\right\vert }%
\ \ \ \ \ \ \ \ \ \ \left(
\begin{array}
[c]{c}%
\text{here, we have renamed the}\\
\text{summation index }I\text{ as }J
\end{array}
\right) \\
&  =\sum_{\substack{J\subseteq G;\\R\subseteq J;\\g\in J}}\left(  -1\right)
^{\left\vert J\right\vert }+\underbrace{\sum_{\substack{J\subseteq
G;\\R\subseteq J;\\g\notin J}}\left(  -1\right)  ^{\left\vert J\right\vert }%
}_{\substack{=-\sum_{\substack{J\subseteq G;\\R\subseteq J;\\g\in J}}\left(
-1\right)  ^{\left\vert J\right\vert }\\\text{(by
(\ref{pf.lem.sol.det.sumdets1.2.short.15}))}}}=\sum_{\substack{J\subseteq
G;\\R\subseteq J;\\g\in J}}\left(  -1\right)  ^{\left\vert J\right\vert }%
-\sum_{\substack{J\subseteq G;\\R\subseteq J;\\g\in J}}\left(  -1\right)
^{\left\vert J\right\vert }=0.
\end{align*}
Comparing this with $\left(  -1\right)  ^{\left\vert G\right\vert
}\underbrace{\left[  G\setminus H\subseteq T\right]  }%
_{\substack{=0\\\text{(since }G\setminus H\not \subseteq T\text{)}}}=0$, we
obtain%
\[
\sum_{\substack{I\subseteq G;\\H\subseteq I;\\T\subseteq I}}\left(  -1\right)
^{\left\vert I\right\vert }=\left(  -1\right)  ^{\left\vert G\right\vert
}\left[  G\setminus H\subseteq T\right]  .
\]
Hence, Lemma \ref{lem.sol.det.sumdets1.2} is proven in Case 1.
\end{vershort}

\begin{verlong}
Let $\mathcal{P}\left(  G\right)  $ be the set of all subsets of $G$. Define
two subsets $A$ and $B$ of $\mathcal{P}\left(  G\right)  $ by%
\begin{align*}
A  &  =\left\{  I\in\mathcal{P}\left(  G\right)  \ \mid\ R\subseteq I\text{
and }g\in I\right\}  \ \ \ \ \ \ \ \ \ \ \text{and}\\
B  &  =\left\{  I\in\mathcal{P}\left(  G\right)  \ \mid\ R\subseteq I\text{
and }g\notin I\right\}  .
\end{align*}

For each $J\in A$, we have $J\setminus\left\{  g\right\}  \in B$%
\ \ \ \ \footnote{\textit{Proof.} Let $J\in A$. Thus, $J\in A=\left\{
I\in\mathcal{P}\left(  G\right)  \ \mid\ R\subseteq I\text{ and }g\in
I\right\}  $. In other words, $J$ is an $I\in\mathcal{P}\left(  G\right)  $
satisfying $R\subseteq I$ and $g\in I$. In other words, $J$ is an element of
$\mathcal{P}\left(  G\right)  $ and satisfies $R\subseteq J$ and $g\in J$.
Notice that $J$ is an element of $\mathcal{P}\left(  G\right)  $; in other
words, $J$ is a subset of $G$ (since $\mathcal{P}\left(  G\right)  $ is the
set of all subsets of $G$). Hence, $J\setminus\left\{  g\right\}  $ is a
subset of $G$ (since $J\setminus\left\{  g\right\}  $ is a subset of $J$). In
other words, $J\setminus\left\{  g\right\}  $ is an element of $\mathcal{P}%
\left(  G\right)  $ (since $\mathcal{P}\left(  G\right)  $ is the set of all
subsets of $G$). Also, $R\setminus\left\{  g\right\}  =R$ (since $g\notin R$),
so that $R=\underbrace{R}_{\subseteq J}\setminus\left\{  g\right\}  \subseteq
J\setminus\left\{  g\right\}  $. Finally, $g\notin J\setminus\left\{
g\right\}  $ (since $g\in\left\{  g\right\}  $). Hence, $J\setminus\left\{
g\right\}  $ is an element of $\mathcal{P}\left(  G\right)  $ and satisfies
$R\subseteq J\setminus\left\{  g\right\}  $ and $g\notin J\setminus\left\{
g\right\}  $. In other words, $J\setminus\left\{  g\right\}  $ is an
$I\in\mathcal{P}\left(  G\right)  $ satisfying $R\subseteq I$ and $g\notin I$.
In other words, $J\setminus\left\{  g\right\}  \in\left\{  I\in\mathcal{P}%
\left(  G\right)  \ \mid\ R\subseteq I\text{ and }g\notin I\right\}  $. In
view of $B=\left\{  I\in\mathcal{P}\left(  G\right)  \ \mid\ R\subseteq
I\text{ and }g\notin I\right\}  $, this rewrites as $J\setminus\left\{
g\right\}  \in B$. Qed.}. Hence, we can define a map $\alpha:A\rightarrow B$
by%
\[
\left(  \alpha\left(  J\right)  =J\setminus\left\{  g\right\}
\ \ \ \ \ \ \ \ \ \ \text{for each }J\in A\right)  .
\]
Consider this $\alpha$.

For each $J\in B$, we have $J\cup\left\{  g\right\}  \in A$%
\ \ \ \ \footnote{\textit{Proof.} Let $J\in B$. Thus, $J\in B=\left\{
I\in\mathcal{P}\left(  G\right)  \ \mid\ R\subseteq I\text{ and }g\notin
I\right\}  $. In other words, $J$ is an $I\in\mathcal{P}\left(  G\right)  $
satisfying $R\subseteq I$ and $g\notin I$. In other words, $J$ is an element
of $\mathcal{P}\left(  G\right)  $ and satisfies $R\subseteq J$ and $g\notin
J$. Notice that $J$ is an element of $\mathcal{P}\left(  G\right)  $; in other
words, $J$ is a subset of $G$ (since $\mathcal{P}\left(  G\right)  $ is the
set of all subsets of $G$). Hence, $J\cup\left\{  g\right\}  $ is a subset of
$G$ (since $J$ and $\left\{  g\right\}  $ are subsets of $G$). In other words,
$J\cup\left\{  g\right\}  $ is an element of $\mathcal{P}\left(  G\right)  $
(since $\mathcal{P}\left(  G\right)  $ is the set of all subsets of $G$).
Also, $R\subseteq J\subseteq J\cup\left\{  g\right\}  $ and $g\in\left\{
g\right\}  \subseteq J\cup\left\{  g\right\}  $. Hence, $J\cup\left\{
g\right\}  $ is an element of $\mathcal{P}\left(  G\right)  $ and satisfies
$R\subseteq J\cup\left\{  g\right\}  $ and $g\in J\cup\left\{  g\right\}  $.
In other words, $J\cup\left\{  g\right\}  $ is an $I\in\mathcal{P}\left(
G\right)  $ satisfying $R\subseteq I$ and $g\in I$. In other words,
$J\cup\left\{  g\right\}  \in\left\{  I\in\mathcal{P}\left(  G\right)
\ \mid\ R\subseteq I\text{ and }g\in I\right\}  $. In view of $A=\left\{
I\in\mathcal{P}\left(  G\right)  \ \mid\ R\subseteq I\text{ and }g\in
I\right\}  $, this rewrites as $J\cup\left\{  g\right\}  \in A$. Qed.}. Hence,
we can define a map $\beta:B\rightarrow A$ by%
\[
\left(  \beta\left(  J\right)  =J\cup\left\{  g\right\}
\ \ \ \ \ \ \ \ \ \ \text{for each }J\in B\right)  .
\]
Consider this $\beta$.

We have
\begin{equation}
\left(  -1\right)  ^{\left\vert \alpha\left(  J\right)  \right\vert }=-\left(
-1\right)  ^{\left\vert J\right\vert }
\label{pf.lem.sol.det.sumdets1.2.signrev}%
\end{equation}
for each $J\in A$\ \ \ \ \footnote{\textit{Proof of
(\ref{pf.lem.sol.det.sumdets1.2.signrev}):} Let $J\in A$. Thus, $J\in
A=\left\{  I\in\mathcal{P}\left(  G\right)  \ \mid\ R\subseteq I\text{ and
}g\in I\right\}  $. In other words, $J$ is an $I\in\mathcal{P}\left(
G\right)  $ satisfying $R\subseteq I$ and $g\in I$. In other words, $J$ is an
element of $\mathcal{P}\left(  G\right)  $ and satisfies $R\subseteq J$ and
$g\in J$. The definition of $\alpha$ yields $\alpha\left(  J\right)
=J\setminus\left\{  g\right\}  $. Hence, $\left\vert \alpha\left(  J\right)
\right\vert =\left\vert J\setminus\left\{  g\right\}  \right\vert =\left\vert
J\right\vert -1$ (since $g\in J$). Thus, $\left(  -1\right)  ^{\left\vert
\alpha\left(  J\right)  \right\vert }=\left(  -1\right)  ^{\left\vert
J\right\vert -1}=-\left(  -1\right)  ^{\left\vert J\right\vert }$. This proves
(\ref{pf.lem.sol.det.sumdets1.2.signrev}).}.

Also, we have $\beta\circ\alpha=\operatorname*{id}$%
\ \ \ \ \footnote{\textit{Proof.} Let $J\in A$. Thus, $J\in A=\left\{
I\in\mathcal{P}\left(  G\right)  \ \mid\ R\subseteq I\text{ and }g\in
I\right\}  $. In other words, $J$ is an $I\in\mathcal{P}\left(  G\right)  $
satisfying $R\subseteq I$ and $g\in I$. In other words, $J$ is an element of
$\mathcal{P}\left(  G\right)  $ and satisfies $R\subseteq J$ and $g\in J$. The
definition of $\alpha$ yields $\alpha\left(  J\right)  =J\setminus\left\{
g\right\}  $. Now,%
\begin{align*}
\left(  \beta\circ\alpha\right)  \left(  J\right)   &  =\beta\left(
\alpha\left(  J\right)  \right)  =\underbrace{\alpha\left(  J\right)
}_{=J\setminus\left\{  g\right\}  }\cup\left\{  g\right\}
\ \ \ \ \ \ \ \ \ \ \left(  \text{by the definition of }\beta\right) \\
&  =\left(  J\setminus\left\{  g\right\}  \right)  \cup\left\{  g\right\}
=J\ \ \ \ \ \ \ \ \ \ \left(  \text{since }g\in J\right) \\
&  =\operatorname*{id}\left(  J\right)  .
\end{align*}
\par
Now, forget that we fixed $J$. We thus have shown that $\left(  \beta
\circ\alpha\right)  \left(  J\right)  =\operatorname*{id}\left(  J\right)  $
for each $J\in A$. In other words, $\beta\circ\alpha=\operatorname*{id}$.} and
$\alpha\circ\beta=\operatorname*{id}$\ \ \ \ \footnote{\textit{Proof.} Let
$J\in B$. Thus, $J\in B=\left\{  I\in\mathcal{P}\left(  G\right)
\ \mid\ R\subseteq I\text{ and }g\notin I\right\}  $. In other words, $J$ is
an $I\in\mathcal{P}\left(  G\right)  $ satisfying $R\subseteq I$ and $g\notin
I$. In other words, $J$ is an element of $\mathcal{P}\left(  G\right)  $ and
satisfies $R\subseteq J$ and $g\notin J$. The definition of $\beta$ yields
$\beta\left(  J\right)  =J\cup\left\{  g\right\}  $. Now,%
\begin{align*}
\left(  \alpha\circ\beta\right)  \left(  J\right)   &  =\alpha\left(
\beta\left(  J\right)  \right)  =\underbrace{\beta\left(  J\right)  }%
_{=J\cup\left\{  g\right\}  }\setminus\left\{  g\right\}
\ \ \ \ \ \ \ \ \ \ \left(  \text{by the definition of }\alpha\right) \\
&  =\left(  J\cup\left\{  g\right\}  \right)  \setminus\left\{  g\right\}
=J\ \ \ \ \ \ \ \ \ \ \left(  \text{since }g\notin J\right) \\
&  =\operatorname*{id}\left(  J\right)  .
\end{align*}
\par
Now, forget that we fixed $J$. We thus have shown that $\left(  \alpha
\circ\beta\right)  \left(  J\right)  =\operatorname*{id}\left(  J\right)  $
for each $J\in B$. In other words, $\alpha\circ\beta=\operatorname*{id}$.}.
Combining these two equalities, we conclude that the maps $\alpha$ and $\beta$
are mutually inverse. Hence, these maps are invertible, i.e., are bijections.
Thus, in particular, $\alpha$ is a bijection. Hence, we can substitute
$\alpha\left(  J\right)  $ for $I$ in the sum $\sum_{I\in B}\left(  -1\right)
^{\left\vert I\right\vert }$. As a result, we find%
\begin{align*}
\sum_{I\in B}\left(  -1\right)  ^{\left\vert I\right\vert }  &  =\sum_{J\in
A}\underbrace{\left(  -1\right)  ^{\left\vert \alpha\left(  J\right)
\right\vert }}_{\substack{=-\left(  -1\right)  ^{\left\vert J\right\vert
}\\\text{(by (\ref{pf.lem.sol.det.sumdets1.2.signrev}))}}}=\sum_{J\in
A}\left(  -\left(  -1\right)  ^{\left\vert J\right\vert }\right)  =\sum_{I\in
A}\left(  -\left(  -1\right)  ^{\left\vert I\right\vert }\right) \\
&  \ \ \ \ \ \ \ \ \ \ \left(  \text{here, we have renamed the summation index
}J\text{ as }I\right) \\
&  =-\sum_{I\in A}\left(  -1\right)  ^{\left\vert I\right\vert }.
\end{align*}
In other words,%
\begin{equation}
\sum_{I\in A}\left(  -1\right)  ^{\left\vert I\right\vert }+\sum_{I\in
B}\left(  -1\right)  ^{\left\vert I\right\vert }=0.
\label{pf.lem.sol.det.sumdets1.2.sum+sum=0}%
\end{equation}

Now,%
\begin{align}
\underbrace{\sum_{\substack{I\subseteq G;\\H\subseteq I;\\T\subseteq I}%
}}_{\substack{=\sum_{\substack{I\subseteq G;\\R\subseteq I}}\\\text{(by
(\ref{pf.lem.sol.det.sumdets1.2.sumeq2}))}}}\left(  -1\right)  ^{\left\vert
I\right\vert }  &  =\underbrace{\sum_{\substack{I\subseteq G;\\R\subseteq I}%
}}_{\substack{=\sum_{\substack{I\in\mathcal{P}\left(  G\right)  ;\\R\subseteq
I}}\\\text{(since the subsets of }G\\\text{are precisely the elements of
}\mathcal{P}\left(  G\right)  \text{)}}}\left(  -1\right)  ^{\left\vert
I\right\vert }=\sum_{\substack{I\in\mathcal{P}\left(  G\right)  ;\\R\subseteq
I}}\left(  -1\right)  ^{\left\vert I\right\vert }\nonumber\\
&  =\underbrace{\sum_{\substack{I\in\mathcal{P}\left(  G\right)  ;\\R\subseteq
I;\\g\in I}}}_{\substack{=\sum_{I\in A}\\\text{(since }A=\left\{
I\in\mathcal{P}\left(  G\right)  \ \mid\ R\subseteq I\text{ and }g\in
I\right\}  \text{)} }} \left(  -1\right)  ^{\left\vert I\right\vert
}+\underbrace{\sum_{\substack{I\in\mathcal{P}\left(  G\right)  ;\\R\subseteq
I;\\g\notin I}}}_{\substack{=\sum_{I\in B}\\\text{(since }B=\left\{
I\in\mathcal{P}\left(  G\right)  \ \mid\ R\subseteq I\text{ and }g\notin
I\right\}  \text{)} }}\left(  -1\right)  ^{\left\vert I\right\vert
}\nonumber\\
&  \ \ \ \ \ \ \ \ \ \ \left(
\begin{array}
[c]{c}%
\text{since each }I\in\mathcal{P}\left(  G\right)  \text{ must satisfy}\\
\text{either }g\in I\text{ or }g\notin I\text{ (but not both)}%
\end{array}
\right) \nonumber\\
&  =\sum_{I\in A}\left(  -1\right)  ^{\left\vert I\right\vert }+\sum_{I\in
B}\left(  -1\right)  ^{\left\vert I\right\vert }=0
\label{pf.lem.sol.det.sumdets1.2.6}%
\end{align}
(by (\ref{pf.lem.sol.det.sumdets1.2.sum+sum=0})).

On the other hand, $G\setminus H\not \subseteq T$. Hence, the statement
$G\setminus H\subseteq T$ is false. Thus, $\left[  G\setminus H\subseteq
T\right]  =0$, so that $\left(  -1\right)  ^{\left\vert G\right\vert
}\underbrace{\left[  G\setminus H\subseteq T\right]  }_{=0}=0$. Comparing this
with (\ref{pf.lem.sol.det.sumdets1.2.6}), we obtain%
\[
\sum_{\substack{I\subseteq G;\\H\subseteq I;\\T\subseteq I}}\left(  -1\right)
^{\left\vert I\right\vert }=\left(  -1\right)  ^{\left\vert G\right\vert
}\left[  G\setminus H\subseteq T\right]  .
\]
Hence, Lemma \ref{lem.sol.det.sumdets1.2} is proven in Case 1.
\end{verlong}

Let us now consider Case 2. In this case, we have $G\setminus H\subseteq T$.
Hence, $T\supseteq G\setminus H$. Now,%
\[
R=H\cup\underbrace{T}_{\supseteq G\setminus H}\supseteq H\cup\left(
G\setminus H\right)  =G
\]
(since $H$ is a subset of $G$). Combining this with $R\subseteq G$, we obtain
$R=G$.

Now, the equality (\ref{pf.lem.sol.det.sumdets1.2.sumeq2}) becomes%
\begin{align*}
\sum_{\substack{I\subseteq G;\\H\subseteq I;\\T\subseteq I}}  &
=\sum_{\substack{I\subseteq G;\\R\subseteq I}}=\sum_{\substack{I\subseteq
G;\\G\subseteq I}}\ \ \ \ \ \ \ \ \ \ \left(  \text{since }R=G\right) \\
&  =\sum_{\substack{I\subseteq G;\\I=G}}\ \ \ \ \ \ \ \ \ \ \left(
\begin{array}
[c]{c}%
\text{because for a subset }I\text{ of }G\text{, the}\\
\text{statement }\left(  G\subseteq I\right)  \text{ is equivalent to }\left(
I=G\right)
\end{array}
\right)  .
\end{align*}
Hence,%
\[
\underbrace{\sum_{\substack{I\subseteq G;\\H\subseteq I;\\T\subseteq I}%
}}_{=\sum_{\substack{I\subseteq G;\\I=G}}}\left(  -1\right)  ^{\left\vert
I\right\vert }=\sum_{\substack{I\subseteq G;\\I=G}}\left(  -1\right)
^{\left\vert I\right\vert }=\left(  -1\right)  ^{\left\vert G\right\vert
}\ \ \ \ \ \ \ \ \ \ \left(  \text{since }G\text{ is a subset of }G\right)  .
\]
Comparing this with $\left(  -1\right)  ^{\left\vert G\right\vert
}\underbrace{\left[  G\setminus H\subseteq T\right]  }%
_{\substack{=1\\\text{(since }G\setminus H\subseteq T\text{)}}}=\left(
-1\right)  ^{\left\vert G\right\vert }$, we obtain%
\[
\sum_{\substack{I\subseteq G;\\H\subseteq I;\\T\subseteq I}}\left(  -1\right)
^{\left\vert I\right\vert }=\left(  -1\right)  ^{\left\vert G\right\vert
}\left[  G\setminus H\subseteq T\right]  .
\]
Hence, Lemma \ref{lem.sol.det.sumdets1.2} is proven in Case 2.

We have thus proven Lemma \ref{lem.sol.det.sumdets1.2} in each of the two
Cases 1 and 2. Since these two Cases cover all possibilities, we thus conclude
that Lemma \ref{lem.sol.det.sumdets1.2} always holds.
\end{proof}

\begin{proof}
[Solution to Exercise \ref{exe.det.sumdets1}.]\textbf{(a)} We have%
\begin{align*}
&  \sum_{\substack{I\subseteq G;\\H\subseteq I}}\left(  -1\right)
^{\left\vert I\right\vert }\underbrace{b_{I,1}b_{I,2}\cdots b_{I,n}}%
_{=\prod_{j=1}^{n}b_{I,j}}\\
&  =\sum_{\substack{I\subseteq G;\\H\subseteq I}}\left(  -1\right)
^{\left\vert I\right\vert }\prod_{j=1}^{n}\underbrace{b_{I,j}}%
_{\substack{=\sum_{i\in I}b_{i,j}\\\text{(by the definition of }%
b_{I,j}\text{)}}}=\sum_{\substack{I\subseteq G;\\H\subseteq I}}\left(
-1\right)  ^{\left\vert I\right\vert }\underbrace{\prod_{j=1}^{n}\sum_{i\in
I}b_{i,j}}_{\substack{=\sum_{\substack{f:\left[  n\right]  \rightarrow
G;\\f\left(  \left[  n\right]  \right)  \subseteq I}}\prod_{j=1}%
^{n}b_{f\left(  j\right)  ,j}\\\text{(by Lemma \ref{lem.sol.det.sumdets1.1})}%
}}\\
&  =\sum_{\substack{I\subseteq G;\\H\subseteq I}}\left(  -1\right)
^{\left\vert I\right\vert }\sum_{\substack{f:\left[  n\right]  \rightarrow
G;\\f\left(  \left[  n\right]  \right)  \subseteq I}}\prod_{j=1}%
^{n}b_{f\left(  j\right)  ,j}=\underbrace{\sum_{\substack{I\subseteq
G;\\H\subseteq I}}\sum_{\substack{f:\left[  n\right]  \rightarrow G;\\f\left(
\left[  n\right]  \right)  \subseteq I}}}_{=\sum_{f:\left[  n\right]
\rightarrow G}\sum_{\substack{I\subseteq G;\\H\subseteq I;\\f\left(  \left[
n\right]  \right)  \subseteq I}}}\left(  -1\right)  ^{\left\vert I\right\vert
}\prod_{j=1}^{n}b_{f\left(  j\right)  ,j}\\
&  =\sum_{f:\left[  n\right]  \rightarrow G}\sum_{\substack{I\subseteq
G;\\H\subseteq I;\\f\left(  \left[  n\right]  \right)  \subseteq I}}\left(
-1\right)  ^{\left\vert I\right\vert }\prod_{j=1}^{n}b_{f\left(  j\right)
,j}=\sum_{f:\left[  n\right]  \rightarrow G}\underbrace{\left(  \sum
_{\substack{I\subseteq G;\\H\subseteq I;\\f\left(  \left[  n\right]  \right)
\subseteq I}}\left(  -1\right)  ^{\left\vert I\right\vert }\right)
}_{\substack{=\left(  -1\right)  ^{\left\vert G\right\vert }\left[  G\setminus
H\subseteq f\left(  \left[  n\right]  \right)  \right]  \\\text{(by Lemma
\ref{lem.sol.det.sumdets1.2}}\\\text{(applied to }T=f\left(  \left[  n\right]
\right)  \text{))}}}\prod_{j=1}^{n}b_{f\left(  j\right)  ,j}\\
&  =\sum_{f:\left[  n\right]  \rightarrow G}\left(  -1\right)  ^{\left\vert
G\right\vert }\left[  G\setminus H\subseteq f\left(  \left[  n\right]
\right)  \right]  \prod_{j=1}^{n}b_{f\left(  j\right)  ,j}\\
&  =\sum_{\substack{f:\left[  n\right]  \rightarrow G;\\G\setminus H\subseteq
f\left(  \left[  n\right]  \right)  }}\left(  -1\right)  ^{\left\vert
G\right\vert }\underbrace{\left[  G\setminus H\subseteq f\left(  \left[
n\right]  \right)  \right]  }_{\substack{=1\\\text{(since }G\setminus
H\subseteq f\left(  \left[  n\right]  \right)  \text{)}}}\prod_{j=1}%
^{n}b_{f\left(  j\right)  ,j}\\
&  \ \ \ \ \ \ \ \ \ \ +\sum_{\substack{f:\left[  n\right]  \rightarrow
G;\\\text{not }G\setminus H\subseteq f\left(  \left[  n\right]  \right)
}}\left(  -1\right)  ^{\left\vert G\right\vert }\underbrace{\left[  G\setminus
H\subseteq f\left(  \left[  n\right]  \right)  \right]  }%
_{\substack{=0\\\text{(since we don't have }G\setminus H\subseteq f\left(
\left[  n\right]  \right)  \text{)}}}\prod_{j=1}^{n}b_{f\left(  j\right)
,j}\\
&  \ \ \ \ \ \ \ \ \ \ \left(
\begin{array}
[c]{c}%
\text{since each map }f:\left[  n\right]  \rightarrow G\text{ satisfies either
}\left(  G\setminus H\subseteq f\left(  \left[  n\right]  \right)  \right) \\
\text{or }\left(  \text{not }G\setminus H\subseteq f\left(  \left[  n\right]
\right)  \right)  \text{ (but not both)}%
\end{array}
\right)
\end{align*}%
\begin{align*}
&  =\sum_{\substack{f:\left[  n\right]  \rightarrow G;\\G\setminus H\subseteq
f\left(  \left[  n\right]  \right)  }}\left(  -1\right)  ^{\left\vert
G\right\vert }\prod_{j=1}^{n}b_{f\left(  j\right)  ,j}+\underbrace{\sum
_{\substack{f:\left[  n\right]  \rightarrow G;\\\text{not }G\setminus
H\subseteq f\left(  \left[  n\right]  \right)  }}\left(  -1\right)
^{\left\vert G\right\vert }0\prod_{j=1}^{n}b_{f\left(  j\right)  ,j}}_{=0}\\
&  =\sum_{\substack{f:\left[  n\right]  \rightarrow G;\\G\setminus H\subseteq
f\left(  \left[  n\right]  \right)  }}\left(  -1\right)  ^{\left\vert
G\right\vert }\prod_{j=1}^{n}b_{f\left(  j\right)  ,j}=\left(  -1\right)
^{\left\vert G\right\vert }\sum_{\substack{f:\left[  n\right]  \rightarrow
G;\\G\setminus H\subseteq f\left(  \left[  n\right]  \right)  }%
}\underbrace{\prod_{j=1}^{n}b_{f\left(  j\right)  ,j}}_{=b_{f\left(  1\right)
,1}b_{f\left(  2\right)  ,2}\cdots b_{f\left(  n\right)  ,n}}\\
&  =\left(  -1\right)  ^{\left\vert G\right\vert }\sum_{\substack{f:\left[
n\right]  \rightarrow G;\\G\setminus H\subseteq f\left(  \left[  n\right]
\right)  }}b_{f\left(  1\right)  ,1}b_{f\left(  2\right)  ,2}\cdots
b_{f\left(  n\right)  ,n}.
\end{align*}
This solves Exercise \ref{exe.det.sumdets1} \textbf{(a)}.

\textbf{(b)} Assume that $n<\left\vert G\setminus H\right\vert $. Then, there
exists no $f:\left[  n\right]  \rightarrow G$ satisfying $G\setminus
H\subseteq f\left(  \left[  n\right]  \right)  $%
\ \ \ \ \footnote{\textit{Proof.} Let $f:\left[  n\right]  \rightarrow G$ be
such that $G\setminus H\subseteq f\left(  \left[  n\right]  \right)  $.\ We
shall derive a contradiction.
\par
Clearly, $\left[  n\right]  $ is a subset of the finite set $\left[  n\right]
$. Thus, Lemma \ref{lem.jectivity.pigeon0} \textbf{(a)} (applied to $\left[
n\right]  $, $G$ and $\left[  n\right]  $ instead of $U$, $V$ and $S$) yields
$\left\vert f\left(  \left[  n\right]  \right)  \right\vert \leq\left\vert
\left[  n\right]  \right\vert =n$. But $G\setminus H\subseteq f\left(  \left[
n\right]  \right)  $; hence, $\left\vert G\setminus H\right\vert
\leq\left\vert f\left(  \left[  n\right]  \right)  \right\vert \leq n$. Now,
$n<\left\vert G\setminus H\right\vert \leq n$. This is absurd. Hence, we have
obtained a contradiction.
\par
Now, forget that we fixed $f$. We thus have found a contradiction for each
$f:\left[  n\right]  \rightarrow G$ satisfying $G\setminus H\subseteq f\left(
\left[  n\right]  \right)  $. Hence, there exists no $f:\left[  n\right]
\rightarrow G$ satisfying $G\setminus H\subseteq f\left(  \left[  n\right]
\right)  $.}. Hence, the sum
\[
\sum_{\substack{f:\left[  n\right]  \rightarrow G;\\G\setminus H\subseteq
f\left(  \left[  n\right]  \right)  }}b_{f\left(  1\right)  ,1}b_{f\left(
2\right)  ,2}\cdots b_{f\left(  n\right)  ,n}%
\]
is empty, and thus equals $0$. In other words,%
\[
\sum_{\substack{f:\left[  n\right]  \rightarrow G;\\G\setminus H\subseteq
f\left(  \left[  n\right]  \right)  }}b_{f\left(  1\right)  ,1}b_{f\left(
2\right)  ,2}\cdots b_{f\left(  n\right)  ,n}=0.
\]
Now, Exercise \ref{exe.det.sumdets1} \textbf{(a)} yields%
\[
\sum_{\substack{I\subseteq G;\\H\subseteq I}}\left(  -1\right)  ^{\left\vert
I\right\vert }b_{I,1}b_{I,2}\cdots b_{I,n}=\left(  -1\right)  ^{\left\vert
G\right\vert }\underbrace{\sum_{\substack{f:\left[  n\right]  \rightarrow
G;\\G\setminus H\subseteq f\left(  \left[  n\right]  \right)  }}b_{f\left(
1\right)  ,1}b_{f\left(  2\right)  ,2}\cdots b_{f\left(  n\right)  ,n}}%
_{=0}=0.
\]
This solves Exercise \ref{exe.det.sumdets1} \textbf{(b)}.

\textbf{(c)} Assume that $n=\left\vert G\right\vert $. For any map $f:\left[
n\right]  \rightarrow G$, we have the following logical equivalence:%
\[
\left(  G\setminus\varnothing\subseteq f\left(  \left[  n\right]  \right)
\right)  \ \Longleftrightarrow\ \left(  f\text{ is bijective}\right)
\]
\footnote{\textit{Proof.} We know that $\left[  n\right]  $ and $G$ are two
finite sets such that $\left\vert \left[  n\right]  \right\vert \leq\left\vert
G\right\vert $ (since $\left\vert \left[  n\right]  \right\vert =n \leq n
=\left\vert G\right\vert $). Thus, Lemma \ref{lem.jectivity.pigeon-surj}
(applied to $U=\left[  n\right]  $ and $V=G$) shows that we have the following
logical equivalence:%
\[
\left(  f\text{ is surjective}\right)  \ \Longleftrightarrow\ \left(  f\text{
is bijective}\right)  .
\]
\par
But we have the following logical equivalence:%
\begin{align*}
\left(  \underbrace{G\setminus\varnothing}_{=G}\subseteq f\left(  \left[
n\right]  \right)  \right)  \  &  \Longleftrightarrow\ \left(  G\subseteq
f\left(  \left[  n\right]  \right)  \right)  \ \Longleftrightarrow\ \left(
G=f\left(  \left[  n\right]  \right)  \right) \\
&  \ \ \ \ \ \ \ \ \ \ \left(  \text{since }f\left(  \left[  n\right]
\right)  \text{ is a subset of }G\right) \\
&  \Longleftrightarrow\ \left(  f\text{ is surjective}\right)
\ \Longleftrightarrow\ \left(  f\text{ is bijective}\right)  .
\end{align*}
Qed.}. Thus, we have the following equality of summation signs:%
\begin{equation}
\sum_{\substack{f:\left[  n\right]  \rightarrow G;\\G\setminus\varnothing
\subseteq f\left(  \left[  n\right]  \right)  }}=\sum_{\substack{f:\left[
n\right]  \rightarrow G;\\f\text{ is bijective}}}.
\label{sol.det.sumdets1.c.1}%
\end{equation}

Also, each subset $I$ of $G$ satisfies $\varnothing\subseteq I$. Hence, we
have the following equality of summation signs:%
\begin{equation}
\sum_{\substack{I\subseteq G;\\\varnothing\subseteq I}}=\sum_{I\subseteq G}.
\label{sol.det.sumdets1.c.2}%
\end{equation}

Now, $\varnothing$ is a subset of $G$. Hence, Exercise \ref{exe.det.sumdets1}
\textbf{(a)} (applied to $\varnothing$ instead of $H$) yields%
\begin{align}
\sum_{\substack{I\subseteq G;\\\varnothing\subseteq I}}\left(  -1\right)
^{\left\vert I\right\vert }b_{I,1}b_{I,2}\cdots b_{I,n}  &  =\left(
-1\right)  ^{\left\vert G\right\vert }\underbrace{\sum_{\substack{f:\left[
n\right]  \rightarrow G;\\G\setminus\varnothing\subseteq f\left(  \left[
n\right]  \right)  }}}_{\substack{=\sum_{\substack{f:\left[  n\right]
\rightarrow G;\\f\text{ is bijective}}}\\\text{(by (\ref{sol.det.sumdets1.c.1}%
))}}}b_{f\left(  1\right)  ,1}b_{f\left(  2\right)  ,2}\cdots b_{f\left(
n\right)  ,n}\nonumber\\
&  =\left(  -1\right)  ^{\left\vert G\right\vert }\sum_{\substack{f:\left[
n\right]  \rightarrow G;\\f\text{ is bijective}}}b_{f\left(  1\right)
,1}b_{f\left(  2\right)  ,2}\cdots b_{f\left(  n\right)  ,n}.
\label{sol.det.sumdets1.c.4}%
\end{align}

Now,%
\begin{align*}
&  \sum_{I\subseteq G}\underbrace{\left(  -1\right)  ^{\left\vert G\setminus
I\right\vert }}_{\substack{=\left(  -1\right)  ^{\left\vert G\right\vert
-\left\vert I\right\vert }\\\text{(since }\left\vert G\setminus I\right\vert
=\left\vert G\right\vert -\left\vert I\right\vert \\\text{(because }I\subseteq
G\text{))}}}b_{I,1}b_{I,2}\cdots b_{I,n}\\
&  =\sum_{I\subseteq G}\underbrace{\left(  -1\right)  ^{\left\vert
G\right\vert -\left\vert I\right\vert }}_{\substack{=\left(  -1\right)
^{\left\vert G\right\vert +\left\vert I\right\vert }\\\text{(since }\left\vert
G\right\vert -\left\vert I\right\vert \equiv\left\vert G\right\vert
+\left\vert I\right\vert \operatorname{mod}2\text{)}}}b_{I,1}b_{I,2}\cdots
b_{I,n}=\underbrace{\sum_{I\subseteq G}}_{\substack{=\sum
_{\substack{I\subseteq G;\\\varnothing\subseteq I}}\\\text{(by
(\ref{sol.det.sumdets1.c.2}))}}}\underbrace{\left(  -1\right)  ^{\left\vert
G\right\vert +\left\vert I\right\vert }}_{=\left(  -1\right)  ^{\left\vert
G\right\vert }\left(  -1\right)  ^{\left\vert I\right\vert }}b_{I,1}%
b_{I,2}\cdots b_{I,n}\\
&  =\sum_{\substack{I\subseteq G;\\\varnothing\subseteq I}}\left(  -1\right)
^{\left\vert G\right\vert }\left(  -1\right)  ^{\left\vert I\right\vert
}b_{I,1}b_{I,2}\cdots b_{I,n}=\left(  -1\right)  ^{\left\vert G\right\vert
}\underbrace{\sum_{\substack{I\subseteq G;\\\varnothing\subseteq I}}\left(
-1\right)  ^{\left\vert I\right\vert }b_{I,1}b_{I,2}\cdots b_{I,n}%
}_{\substack{=\left(  -1\right)  ^{\left\vert G\right\vert }\sum
_{\substack{f:\left[  n\right]  \rightarrow G;\\f\text{ is bijective}%
}}b_{f\left(  1\right)  ,1}b_{f\left(  2\right)  ,2}\cdots b_{f\left(
n\right)  ,n}\\\text{(by (\ref{sol.det.sumdets1.c.4}))}}}\\
&  =\underbrace{\left(  -1\right)  ^{\left\vert G\right\vert }\left(
-1\right)  ^{\left\vert G\right\vert }}_{\substack{=\left(  \left(  -1\right)
\left(  -1\right)  \right)  ^{\left\vert G\right\vert }=1^{\left\vert
G\right\vert }\\\text{(since }\left(  -1\right)  \left(  -1\right)
=1\text{)}}}\sum_{\substack{f:\left[  n\right]  \rightarrow G;\\f\text{ is
bijective}}}b_{f\left(  1\right)  ,1}b_{f\left(  2\right)  ,2}\cdots
b_{f\left(  n\right)  ,n}\\
&  =\underbrace{1^{\left\vert G\right\vert }}_{=1}\sum_{\substack{f:\left[
n\right]  \rightarrow G;\\f\text{ is bijective}}}b_{f\left(  1\right)
,1}b_{f\left(  2\right)  ,2}\cdots b_{f\left(  n\right)  ,n}\\
&  =\sum_{\substack{f:\left[  n\right]  \rightarrow G;\\f\text{ is bijective}%
}}b_{f\left(  1\right)  ,1}b_{f\left(  2\right)  ,2}\cdots b_{f\left(
n\right)  ,n}.
\end{align*}
This solves Exercise \ref{exe.det.sumdets1} \textbf{(c)}.
\end{proof}

\subsection{Solution to Exercise \ref{exe.det.sumdets2}}

In this section, we shall use the same notations that we have introduced in
Section \ref{sect.sol.noncomm.polarization}.

Let us begin by stating a simple corollary from Exercise
\ref{exe.det.sumdets1} \textbf{(b)}:

\begin{lemma}
\label{lem.sol.det.sumdets2.1}Let $G$ be a finite set. Let $n\in\mathbb{N}$ be
such that $n<\left\vert G\right\vert $.

For each $i\in G$ and $j\in\left[  n\right]  $, let $b_{i,j}$ be an element of
$\mathbb{K}$. Then,
\[
\sum_{I\subseteq G}\left(  -1\right)  ^{\left\vert I\right\vert }\prod
_{i=1}^{n}\sum_{g\in I}b_{g,i}=0.
\]

\end{lemma}

\begin{proof}
[Proof of Lemma \ref{lem.sol.det.sumdets2.1}.]Clearly, $\varnothing$ is a
subset of $G$. Also, $G\setminus\varnothing=G$, so that $\left\vert
G\setminus\varnothing\right\vert =\left\vert G\right\vert $. Thus,
$n<\left\vert G\right\vert =\left\vert G\setminus\varnothing\right\vert $.

For each $j\in\left[  n\right]  $ and each subset $I$ of $G$, we define an
element $b_{I,j}\in\mathbb{K}$ by $b_{I,j}=\sum_{i\in I}b_{i,j}$. Thus, for
each $j\in\left[  n\right]  $ and each subset $I$ of $G$, we have%
\begin{equation}
b_{I,j}=\sum_{i\in I}b_{i,j}=\sum_{g\in I}b_{g,j}
\label{pf.lem.sol.det.sumdets2.1.1}%
\end{equation}
(here, we have substituted $g$ for $i$ in the sum).

Now, Exercise \ref{exe.det.sumdets1} \textbf{(b)} (applied to $\mathbb{K}$ and
$\varnothing$ instead of $\mathbb{L}$ and $H$) yields%
\begin{equation}
\sum_{\substack{I\subseteq G;\\\varnothing\subseteq I}}\left(  -1\right)
^{\left\vert I\right\vert }b_{I,1}b_{I,2}\cdots b_{I,n}=0.
\label{pf.lem.sol.det.sumdets2.1.2}%
\end{equation}
But each subset $I$ of $G$ satisfies $\varnothing\subseteq I$. Hence, we have
the following equality of summation signs:%
\[
\sum_{\substack{I\subseteq G;\\\varnothing\subseteq I}}=\sum_{I\subseteq G}.
\]
Thus, the equality (\ref{pf.lem.sol.det.sumdets2.1.2}) rewrites as follows:%
\[
\sum_{I\subseteq G}\left(  -1\right)  ^{\left\vert I\right\vert }%
b_{I,1}b_{I,2}\cdots b_{I,n}=0.
\]
Hence,%
\begin{align*}
0  &  =\sum_{I\subseteq G}\left(  -1\right)  ^{\left\vert I\right\vert
}\underbrace{b_{I,1}b_{I,2}\cdots b_{I,n}}_{=\prod_{i=1}^{n}b_{I,i}}%
=\sum_{I\subseteq G}\left(  -1\right)  ^{\left\vert I\right\vert }\prod
_{i=1}^{n}\underbrace{b_{I,i}}_{\substack{=\sum_{g\in I}b_{g,i}\\\text{(by
(\ref{pf.lem.sol.det.sumdets2.1.1}) (applied to }j=i\text{))}}}\\
&  =\sum_{I\subseteq G}\left(  -1\right)  ^{\left\vert I\right\vert }%
\prod_{i=1}^{n}\sum_{g\in I}b_{g,i}.
\end{align*}
This proves Lemma \ref{lem.sol.det.sumdets2.1}.
\end{proof}

\begin{proof}
[Solution to Exercise \ref{exe.det.sumdets2}.]For each $g\in G$, let us write
the $n\times n$-matrix $A_{g}$ in the form
\begin{equation}
A_{g}=\left(  a_{g,i,j}\right)  _{1\leq i\leq n,\ 1\leq j\leq n}.
\label{sol.det.sumdets2.1}%
\end{equation}

Let $I$ be any subset of $G$. Then,%
\[
\sum_{g\in I}\underbrace{A_{g}}_{\substack{=\left(  a_{g,i,j}\right)  _{1\leq
i\leq n,\ 1\leq j\leq n}\\\text{(by (\ref{sol.det.sumdets2.1}))}}}=\sum_{g\in
I}\left(  a_{g,i,j}\right)  _{1\leq i\leq n,\ 1\leq j\leq n}=\left(
\sum_{g\in I}a_{g,i,j}\right)  _{1\leq i\leq n,\ 1\leq j\leq n}.
\]
Hence,%
\[
\det\left(  \sum_{g\in I}A_{g}\right)  =\sum_{\sigma\in S_{n}}\left(
-1\right)  ^{\sigma}\prod_{i=1}^{n}\sum_{g\in I}a_{g,i,\sigma\left(  i\right)
}%
\]
(by (\ref{eq.det.eq.2}) (applied to $\sum_{g\in I}A_{g}$ and $\sum_{g\in
I}a_{g,i,j}$ instead of $A$ and $a_{i,j}$)). In view of%
\[
\sum_{g\in I}A_{g}=\sum_{i\in I}A_{i}\ \ \ \ \ \ \ \ \ \ \left(  \text{here,
we have renamed the summation index }g\text{ as }i\right)  ,
\]
this rewrites as follows:%
\begin{equation}
\det\left(  \sum_{i\in I}A_{i}\right)  =\sum_{\sigma\in S_{n}}\left(
-1\right)  ^{\sigma}\prod_{i=1}^{n}\sum_{g\in I}a_{g,i,\sigma\left(  i\right)
}. \label{sol.det.sumdets2.2}%
\end{equation}

Now, forget that we fixed $I$. We thus have proven (\ref{sol.det.sumdets2.2})
for each subset $I$ of $G$. Now,%
\begin{align*}
&  \sum_{I\subseteq G}\left(  -1\right)  ^{\left\vert I\right\vert
}\underbrace{\det\left(  \sum_{i\in I}A_{i}\right)  }_{\substack{=\sum
_{\sigma\in S_{n}}\left(  -1\right)  ^{\sigma}\prod_{i=1}^{n}\sum_{g\in
I}a_{g,i,\sigma\left(  i\right)  }\\\text{(by (\ref{sol.det.sumdets2.2}))}}}\\
&  =\sum_{I\subseteq G}\left(  -1\right)  ^{\left\vert I\right\vert }%
\sum_{\sigma\in S_{n}}\left(  -1\right)  ^{\sigma}\prod_{i=1}^{n}\sum_{g\in
I}a_{g,i,\sigma\left(  i\right)  }=\underbrace{\sum_{I\subseteq G}\sum
_{\sigma\in S_{n}}}_{=\sum_{\sigma\in S_{n}}\sum_{I\subseteq G}}\left(
-1\right)  ^{\left\vert I\right\vert }\left(  -1\right)  ^{\sigma}\prod
_{i=1}^{n}\sum_{g\in I}a_{g,i,\sigma\left(  i\right)  }\\
&  =\sum_{\sigma\in S_{n}}\sum_{I\subseteq G}\left(  -1\right)  ^{\left\vert
I\right\vert }\left(  -1\right)  ^{\sigma}\prod_{i=1}^{n}\sum_{g\in
I}a_{g,i,\sigma\left(  i\right)  }=\sum_{\sigma\in S_{n}}\left(  -1\right)
^{\sigma}\underbrace{\sum_{I\subseteq G}\left(  -1\right)  ^{\left\vert
I\right\vert }\prod_{i=1}^{n}\sum_{g\in I}a_{g,i,\sigma\left(  i\right)  }%
}_{\substack{=0\\\text{(by Lemma \ref{lem.sol.det.sumdets2.1}}\\\text{(applied
to }b_{i,j}=a_{i,j,\sigma\left(  j\right)  }\text{))}}}\\
&  =\sum_{\sigma\in S_{n}}\left(  -1\right)  ^{\sigma}0=0.
\end{align*}
This solves Exercise \ref{exe.det.sumdets2}.
\end{proof}

\subsection{Solution to Exercise \ref{exe.powerdet.gen}}

We shall first prove some identities in preparation for the solution to
Exercise \ref{exe.powerdet.gen}.

First, we introduce a notation: For every $n\in\mathbb{N}$, let $\left[
n\right]  $ denote the set $\left\{  1,2,\ldots,n\right\}  $. We shall also
use the notation introduced in Definition \ref{def.submatrix.minor} throughout
this section.

Let us now state a simple corollary of Lemma \ref{lem.prodrule2}:

\begin{corollary}
\label{cor.sol.powerdet.gen.prodrule}Let $n\in\mathbb{N}$. Let $A=\left(
a_{i,j}\right)  _{1\leq i\leq n,\ 1\leq j\leq n}$ be an $n\times n$-matrix.
Let $\sigma\in S_{n}$. Then,%
\[
\left(  \sum_{i=1}^{n}a_{i,\sigma\left(  i\right)  }\right)  ^{k}=\sum
_{\kappa:\left[  k\right]  \rightarrow\left[  n\right]  }\prod_{i=1}%
^{k}a_{\kappa\left(  i\right)  ,\sigma\left(  \kappa\left(  i\right)  \right)
}%
\]
for every $k\in\mathbb{N}$.
\end{corollary}

\begin{proof}
[Proof of Corollary \ref{cor.sol.powerdet.gen.prodrule}.]Let $g\in\mathbb{N}$.
Then, Lemma \ref{lem.prodrule2} (applied to $g$, $n$ and $a_{k,\sigma\left(
k\right)  }$ instead of $n$, $m$ and $p_{i,k}$) yields%
\[
\prod_{i=1}^{g}\sum_{k=1}^{n}a_{k,\sigma\left(  k\right)  }=\sum
_{\kappa:\left[  g\right]  \rightarrow\left[  n\right]  }\prod_{i=1}%
^{g}a_{\kappa\left(  i\right)  ,\sigma\left(  \kappa\left(  i\right)  \right)
}.
\]
Thus,%
\begin{align}
\sum_{\kappa:\left[  g\right]  \rightarrow\left[  n\right]  }\prod_{i=1}%
^{g}a_{\kappa\left(  i\right)  ,\sigma\left(  \kappa\left(  i\right)  \right)
}  &  =\prod_{i=1}^{g}\sum_{k=1}^{n}a_{k,\sigma\left(  k\right)  }=\left(
\sum_{k=1}^{n}a_{k,\sigma\left(  k\right)  }\right)  ^{g}\nonumber\\
&  =\left(  \sum_{i=1}^{n}a_{i,\sigma\left(  i\right)  }\right)  ^{g}
\label{pf.cor.sol.powerdet.gen.prodrule.1}%
\end{align}
(here, we have renamed the summation index $k$ as $i$ in the sum).

Now, forget that we fixed $g$. We thus have proven the identity
(\ref{pf.cor.sol.powerdet.gen.prodrule.1}) for each $g\in\mathbb{N}$.

Now, fix $k\in\mathbb{N}$. Then, (\ref{pf.cor.sol.powerdet.gen.prodrule.1})
(applied to $g=k$) yields $\sum_{\kappa:\left[  k\right]  \rightarrow\left[
n\right]  }\prod_{i=1}^{k}a_{\kappa\left(  i\right)  ,\sigma\left(
\kappa\left(  i\right)  \right)  }=\left(  \sum_{i=1}^{n}a_{i,\sigma\left(
i\right)  }\right)  ^{k}$. This proves Corollary
\ref{cor.sol.powerdet.gen.prodrule}.
\end{proof}

Next, we state a simple lemma:

\begin{lemma}
\label{lem.sol.powerdet.gen.toofew}Let $n\in\mathbb{N}$ and $k\in\mathbb{N}$.
Let $\kappa:\left[  k\right]  \rightarrow\left[  n\right]  $ be a map such
that $\left\vert \kappa\left(  \left[  k\right]  \right)  \right\vert <n-1$.
Then,%
\[
\sum_{\sigma\in S_{n}}\left(  -1\right)  ^{\sigma}\prod_{i=1}^{k}%
a_{\kappa\left(  i\right)  ,\sigma\left(  \kappa\left(  i\right)  \right)
}=0.
\]

\end{lemma}

\begin{proof}
[Proof of Lemma \ref{lem.sol.powerdet.gen.toofew}.]We have $\left\vert \left[
n\right]  \setminus\kappa\left(  \left[  k\right]  \right)  \right\vert \geq
2$\ \ \ \ \footnote{\textit{Proof.} We have $\left\vert \kappa\left(  \left[
k\right]  \right)  \right\vert <n-1$. Thus, $\left\vert \kappa\left(  \left[
k\right]  \right)  \right\vert \leq\left(  n-1\right)  -1$ (since both
$\left\vert \kappa\left(  \left[  k\right]  \right)  \right\vert $ and $n-1$
are integers). Hence, $\left\vert \kappa\left(  \left[  k\right]  \right)
\right\vert \leq\left(  n-1\right)  -1=n-2$.
\par
But $\kappa\left(  \left[  k\right]  \right)  \subseteq\left[  n\right]  $.
Hence, $\left\vert \left[  n\right]  \setminus\kappa\left(  \left[  k\right]
\right)  \right\vert =\underbrace{\left\vert \left[  n\right]  \right\vert
}_{=n}-\underbrace{\left\vert \kappa\left(  \left[  k\right]  \right)
\right\vert }_{\leq n-2}\geq n-\left(  n-2\right)  =2$. Qed.}. Hence, the set
$\left[  n\right]  \setminus\kappa\left(  \left[  k\right]  \right)  $ has at
least $2$ elements. In other words, the set $\left[  n\right]  \setminus
\kappa\left(  \left[  k\right]  \right)  $ contains two distinct elements.
Choose two such elements, and denote them by $a$ and $b$. Thus, $a$ and $b$
are two distinct elements of the set $\left[  n\right]  \setminus\kappa\left(
\left[  k\right]  \right)  $.

The elements $a$ and $b$ both belong to $\left[  n\right]  $ (since
$a\in\left[  n\right]  \setminus\kappa\left(  \left[  k\right]  \right)
\subseteq\left[  n\right]  $ and $b\in\left[  n\right]  \setminus\kappa\left(
\left[  k\right]  \right)  \subseteq\left[  n\right]  $).

\begin{vershort}
Thus, $a$ and $b$ are two distinct elements of $\left[  n\right]  =\left\{
1,2,\ldots,n\right\}  $. Hence, a transposition $t_{a,b}\in S_{n}$ is defined
(according to Definition \ref{def.transpos}). This transposition satisfies
$t_{a,b}\circ\kappa=\kappa$\ \ \ \ \footnote{\textit{Proof.} We are going to
show that every $i\in\left[  k\right]  $ satisfies $\left(  t_{a,b}\circ
\kappa\right)  \left(  i\right)  =\kappa\left(  i\right)  $.
\par
So let $i\in\left[  k\right]  $. We shall show that $\left(  t_{a,b}%
\circ\kappa\right)  \left(  i\right)  =\kappa\left(  i\right)  $.
\par
We have $t_{a,b}\left(  j\right)  =j$ for every $j\in\left[  n\right]
\setminus\left\{  a,b\right\}  $ (by the definition of $t_{a,b}$).
\par
But $a$ and $b$ are two elements of the set $\left[  n\right]  \setminus
\kappa\left(  \left[  k\right]  \right)  $. Hence, $\left\{  a,b\right\}
\subseteq\left[  n\right]  \setminus\kappa\left(  \left[  k\right]  \right)
$. If we had $\kappa\left(  i\right)  \in\left\{  a,b\right\}  $, then we
would thus have $\kappa\left(  i\right)  \in\left\{  a,b\right\}
\subseteq\left[  n\right]  \setminus\kappa\left(  \left[  k\right]  \right)
$, which would contradict $\kappa\left(  i\right)  \in\kappa\left(  \left[
k\right]  \right)  $. Hence, we cannot have $\kappa\left(  i\right)
\in\left\{  a,b\right\}  $. Thus, we have $\kappa\left(  i\right)
\notin\left\{  a,b\right\}  $.
\par
Combining $\kappa\left(  i\right)  \in\left[  n\right]  $ with $\kappa\left(
i\right)  \notin\left\{  a,b\right\}  $, we obtain $\kappa\left(  i\right)
\in\left[  n\right]  \setminus\left\{  a,b\right\}  $. But recall that
$t_{a,b}\left(  j\right)  =j$ for every $j\in\left[  n\right]  \setminus
\left\{  a,b\right\}  $. Applying this to $j=\kappa\left(  i\right)  $, we
obtain $t_{a,b}\left(  \kappa\left(  i\right)  \right)  =\kappa\left(
i\right)  $ (since $\kappa\left(  i\right)  \in\left[  n\right]
\setminus\left\{  a,b\right\}  $). Hence, $\left(  t_{a,b}\circ\kappa\right)
\left(  i\right)  =t_{a,b}\left(  \kappa\left(  i\right)  \right)
=\kappa\left(  i\right)  $.
\par
Now, let us forget that we fixed $i$. We thus have shown that $\left(
t_{a,b}\circ\kappa\right)  \left(  i\right)  =\kappa\left(  i\right)  $ for
every $i\in\left[  k\right]  $. In other words, $t_{a,b}\circ\kappa=\kappa$,
qed.}.
\end{vershort}

\begin{verlong}
Thus, $a$ and $b$ are two distinct elements of $\left[  n\right]  =\left\{
1,2,\ldots,n\right\}  $. Hence, a transposition $t_{a,b}\in S_{n}$ is defined
(according to Definition \ref{def.transpos}). This transposition satisfies
$t_{a,b}\circ\kappa=\kappa$\ \ \ \ \footnote{\textit{Proof.} We are going to
show that every $i\in\left[  k\right]  $ satisfies $\left(  t_{a,b}\circ
\kappa\right)  \left(  i\right)  =\kappa\left(  i\right)  $.
\par
So let $i\in\left[  k\right]  $. We shall show that $\left(  t_{a,b}%
\circ\kappa\right)  \left(  i\right)  =\kappa\left(  i\right)  $.
\par
The definition of $t_{a,b}$ shows that $t_{a,b}$ is the permutation in $S_{n}$
which swaps $a$ with $b$ while leaving all other elements of $\left\{
1,2,\ldots,n\right\}  $ unchanged. In other words, we have $t_{a,b}\left(
a\right)  =b$, and $t_{a,b}\left(  b\right)  =a$, and $t_{a,b}\left(
j\right)  =j$ for every $j\in\left[  n\right]  \setminus\left\{  a,b\right\}
$.
\par
Let us first assume (for the sake of contradiction) that $\kappa\left(
i\right)  \in\left\{  a,b\right\}  $.
\par
But $a$ and $b$ are two elements of the set $\left[  n\right]  \setminus
\kappa\left(  \left[  k\right]  \right)  $. Hence, $\left\{  a,b\right\}
\subseteq\left[  n\right]  \setminus\kappa\left(  \left[  k\right]  \right)
$. Thus, $\kappa\left(  i\right)  \in\left\{  a,b\right\}  \subseteq\left[
n\right]  \setminus\kappa\left(  \left[  k\right]  \right)  $. In other words,
$\kappa\left(  i\right)  \in\left[  n\right]  $ and $\kappa\left(  i\right)
\notin\kappa\left(  \left[  k\right]  \right)  $. But $\kappa\left(  i\right)
\notin\kappa\left(  \left[  k\right]  \right)  $ clearly contradicts
$\kappa\left(  \underbrace{i}_{\in\left[  k\right]  }\right)  \in\kappa\left(
\left[  k\right]  \right)  $. Thus, we have obtained a contradiction. Hence,
our assumption (that $\kappa\left(  i\right)  \in\left\{  a,b\right\}  $) was
false. Thus, we have $\kappa\left(  i\right)  \notin\left\{  a,b\right\}  $.
\par
Combining $\kappa\left(  i\right)  \in\left[  n\right]  $ with $\kappa\left(
i\right)  \notin\left\{  a,b\right\}  $, we obtain $\kappa\left(  i\right)
\in\left[  n\right]  \setminus\left\{  a,b\right\}  $. But recall that
$t_{a,b}\left(  j\right)  =j$ for every $j\in\left[  n\right]  \setminus
\left\{  a,b\right\}  $. Applying this to $j=\kappa\left(  i\right)  $, we
obtain $t_{a,b}\left(  \kappa\left(  i\right)  \right)  =\kappa\left(
i\right)  $ (since $\kappa\left(  i\right)  \in\left[  n\right]
\setminus\left\{  a,b\right\}  $). Hence, $\left(  t_{a,b}\circ\kappa\right)
\left(  i\right)  =t_{a,b}\left(  \kappa\left(  i\right)  \right)
=\kappa\left(  i\right)  $.
\par
Now, let us forget that we fixed $i$. We thus have shown that $\left(
t_{a,b}\circ\kappa\right)  \left(  i\right)  =\kappa\left(  i\right)  $ for
every $i\in\left[  k\right]  $. In other words, $t_{a,b}\circ\kappa=\kappa$,
qed.}.
\end{verlong}

Let $A_{n}$ be the set of all even permutations in $S_{n}$. Let $C_{n}$ be the
set of all odd permutations in $S_{n}$.

We have $\sigma\circ t_{a,b}\in C_{n}$ for every $\sigma\in A_{n}%
$\ \ \ \ \footnote{We have already proven this during our proof of Lemma
\ref{lem.det.sigma} \textbf{(b)}.}. Hence, we can define a map $\Phi
:A_{n}\rightarrow C_{n}$ by%
\[
\Phi\left(  \sigma\right)  =\sigma\circ t_{a,b}\ \ \ \ \ \ \ \ \ \ \text{for
every }\sigma\in A_{n}.
\]
Consider this map $\Phi$. Furthermore, we have $\sigma\circ\left(
t_{a,b}\right)  ^{-1}\in A_{n}$ for every $\sigma\in C_{n}$%
\ \ \ \ \footnote{We have already proven this during our proof of Lemma
\ref{lem.det.sigma} \textbf{(b)}.}. Thus, we can define a map $\Psi
:C_{n}\rightarrow A_{n}$ by%
\[
\Psi\left(  \sigma\right)  =\sigma\circ\left(  t_{a,b}\right)  ^{-1}%
\ \ \ \ \ \ \ \ \ \ \text{for every }\sigma\in C_{n}.
\]
Consider this map $\Psi$.

The maps $\Phi$ and $\Psi$ are mutually inverse\footnote{We have already
proven this during our proof of Lemma \ref{lem.det.sigma} \textbf{(b)}.}.
Hence, the map $\Phi$ is a bijection. Moreover, every $\sigma\in A_{n}$ and
$i\in\left\{  1,2,\ldots,k\right\}  $ satisfy%
\begin{equation}
a_{\kappa\left(  i\right)  ,\left(  \Phi\left(  \sigma\right)  \right)
\left(  \kappa\left(  i\right)  \right)  }=a_{\kappa\left(  i\right)
,\sigma\left(  \kappa\left(  i\right)  \right)  }
\label{pf.lem.sol.powerdet.gen.toofew.1}%
\end{equation}
\footnote{\textit{Proof of (\ref{pf.lem.sol.powerdet.gen.toofew.1}):} Let
$\sigma\in A_{n}$ and $i\in\left\{  1,2,\ldots,k\right\}  $. Thus,
$i\in\left\{  1,2,\ldots,k\right\}  =\left[  k\right]  $. Now,%
\[
\left(  \underbrace{\Phi\left(  \sigma\right)  }_{=\sigma\circ t_{a,b}%
}\right)  \left(  \kappa\left(  i\right)  \right)  =\left(  \sigma\circ
t_{a,b}\right)  \left(  \kappa\left(  i\right)  \right)  =\sigma\left(
t_{a,b}\left(  \kappa\left(  i\right)  \right)  \right)  =\sigma\left(
\underbrace{\left(  t_{a,b}\circ\kappa\right)  }_{=\kappa}\left(  i\right)
\right)  =\sigma\left(  \kappa\left(  i\right)  \right)  .
\]
Hence, $a_{\kappa\left(  i\right)  ,\left(  \Phi\left(  \sigma\right)
\right)  \left(  \kappa\left(  i\right)  \right)  }=a_{\kappa\left(  i\right)
,\sigma\left(  \kappa\left(  i\right)  \right)  }$. This proves
(\ref{pf.lem.sol.powerdet.gen.toofew.1}).}.

Now,%
\begin{align*}
&  \sum_{\sigma\in S_{n}}\left(  -1\right)  ^{\sigma}\prod_{i=1}^{k}%
a_{\kappa\left(  i\right)  ,\sigma\left(  \kappa\left(  i\right)  \right)  }\\
&  =\underbrace{\sum_{\substack{\sigma\in S_{n};\\\sigma\text{ is even}}%
}}_{\substack{=\sum_{\sigma\in A_{n}}\\\text{(since }A_{n}\text{ is
the}\\\text{set of all even}\\\text{permutations}\\\text{in }S_{n}\text{)}%
}}\underbrace{\left(  -1\right)  ^{\sigma}}_{\substack{=1\\\text{(since
}\sigma\text{ is even)}}}\prod_{i=1}^{k}a_{\kappa\left(  i\right)
,\sigma\left(  \kappa\left(  i\right)  \right)  }+\underbrace{\sum
_{\substack{\sigma\in S_{n};\\\sigma\text{ is odd}}}}_{\substack{=\sum
_{\sigma\in C_{n}}\\\text{(since }C_{n}\text{ is the}\\\text{set of all
odd}\\\text{permutations}\\\text{in }S_{n}\text{)}}}\underbrace{\left(
-1\right)  ^{\sigma}}_{\substack{=-1\\\text{(since }\sigma\text{ is odd)}%
}}\prod_{i=1}^{k}a_{\kappa\left(  i\right)  ,\sigma\left(  \kappa\left(
i\right)  \right)  }\\
&  \ \ \ \ \ \ \ \ \ \ \left(  \text{since every permutation }\sigma\in
S_{n}\text{ is either even or odd, but not both}\right) \\
&  =\sum_{\sigma\in A_{n}}\prod_{i=1}^{k}a_{\kappa\left(  i\right)
,\sigma\left(  \kappa\left(  i\right)  \right)  }+\sum_{\sigma\in C_{n}%
}\left(  -1\right)  \prod_{i=1}^{k}a_{\kappa\left(  i\right)  ,\sigma\left(
\kappa\left(  i\right)  \right)  }\\
&  =\sum_{\sigma\in A_{n}}\prod_{i=1}^{k}a_{\kappa\left(  i\right)
,\sigma\left(  \kappa\left(  i\right)  \right)  }-\sum_{\sigma\in C_{n}}%
\prod_{i=1}^{k}a_{\kappa\left(  i\right)  ,\sigma\left(  \kappa\left(
i\right)  \right)  }=0,
\end{align*}
since%
\begin{align*}
\sum_{\sigma\in C_{n}}\prod_{i=1}^{k}a_{\kappa\left(  i\right)  ,\sigma\left(
\kappa\left(  i\right)  \right)  }  &  =\sum_{\sigma\in A_{n}}\prod_{i=1}%
^{k}\underbrace{a_{\kappa\left(  i\right)  ,\left(  \Phi\left(  \sigma\right)
\right)  \left(  \kappa\left(  i\right)  \right)  }}_{\substack{=a_{\kappa
\left(  i\right)  ,\sigma\left(  \kappa\left(  i\right)  \right)  }\\\text{(by
(\ref{pf.lem.sol.powerdet.gen.toofew.1}))}}}\\
&  \ \ \ \ \ \ \ \ \ \ \left(
\begin{array}
[c]{c}%
\text{here, we have substituted }\Phi\left(  \sigma\right)  \text{ for }%
\sigma\text{ in}\\
\text{the sum, since the map }\Phi\text{ is a bijection}%
\end{array}
\right) \\
&  =\sum_{\sigma\in A_{n}}\prod_{i=1}^{k}a_{\kappa\left(  i\right)
,\sigma\left(  \kappa\left(  i\right)  \right)  }.
\end{align*}
This proves Lemma \ref{lem.sol.powerdet.gen.toofew}.
\end{proof}

We can now solve part \textbf{(a)} of Exercise \ref{exe.powerdet.gen}:

\begin{proposition}
\label{prop.sol.powerdet.gen.a}Let $n\in\mathbb{N}$. Let $A=\left(
a_{i,j}\right)  _{1\leq i\leq n,\ 1\leq j\leq n}$ be an $n\times n$-matrix.
Then,%
\[
\sum_{\sigma\in S_{n}}\left(  -1\right)  ^{\sigma}\left(  \sum_{i=1}%
^{n}a_{i,\sigma\left(  i\right)  }\right)  ^{k}=0
\]
for each $k\in\left\{  0,1,\ldots,n-2\right\}  $.
\end{proposition}

\begin{proof}
[Proof of Proposition \ref{prop.sol.powerdet.gen.a}.]Let $k\in\left\{
0,1,\ldots,n-2\right\}  $. Then, each each map $\kappa:\left[  k\right]
\rightarrow\left[  n\right]  $ satisfies%
\begin{equation}
\sum_{\sigma\in S_{n}}\left(  -1\right)  ^{\sigma}\prod_{i=1}^{k}%
a_{\kappa\left(  i\right)  ,\sigma\left(  \kappa\left(  i\right)  \right)  }=0
\label{pf.prop.sol.powerdet.gen.a.1}%
\end{equation}
\footnote{\textit{Proof of (\ref{pf.prop.sol.powerdet.gen.a.1}):} Let
$\kappa:\left[  k\right]  \rightarrow\left[  n\right]  $ be a map. Then, Lemma
\ref{lem.jectivity.pigeon0} \textbf{(a)} (applied to $U=\left[  k\right]  $,
$V=\left[  n\right]  $, $f=\kappa$ and $S=\left[  k\right]  $) yields
$\left\vert \kappa\left(  \left[  k\right]  \right)  \right\vert
\leq\left\vert \left[  k\right]  \right\vert =k\leq n-2$ (since $k\in\left\{
0,1,\ldots,n-2\right\}  $), so that $\left\vert \kappa\left(  \left[
k\right]  \right)  \right\vert \leq n-2<n-1$. Hence, Lemma
\ref{lem.sol.powerdet.gen.toofew} yields $\sum_{\sigma\in S_{n}}\left(
-1\right)  ^{\sigma}\prod_{i=1}^{k}a_{\kappa\left(  i\right)  ,\sigma\left(
\kappa\left(  i\right)  \right)  }=0$. This proves
(\ref{pf.prop.sol.powerdet.gen.a.1}).}.

Now,%
\begin{align*}
&  \sum_{\sigma\in S_{n}}\left(  -1\right)  ^{\sigma}\underbrace{\left(
\sum_{i=1}^{n}a_{i,\sigma\left(  i\right)  }\right)  ^{k}}_{\substack{=\sum
_{\kappa:\left[  k\right]  \rightarrow\left[  n\right]  }\prod_{i=1}%
^{k}a_{\kappa\left(  i\right)  ,\sigma\left(  \kappa\left(  i\right)  \right)
}\\\text{(by Corollary \ref{cor.sol.powerdet.gen.prodrule})}}}\\
&  =\sum_{\sigma\in S_{n}}\left(  -1\right)  ^{\sigma}\sum_{\kappa:\left[
k\right]  \rightarrow\left[  n\right]  }\prod_{i=1}^{k}a_{\kappa\left(
i\right)  ,\sigma\left(  \kappa\left(  i\right)  \right)  }=\underbrace{\sum
_{\sigma\in S_{n}}\sum_{\kappa:\left[  k\right]  \rightarrow\left[  n\right]
}}_{=\sum_{\kappa:\left[  k\right]  \rightarrow\left[  n\right]  }\sum
_{\sigma\in S_{n}}}\left(  -1\right)  ^{\sigma}\prod_{i=1}^{k}a_{\kappa\left(
i\right)  ,\sigma\left(  \kappa\left(  i\right)  \right)  }\\
&  =\sum_{\kappa:\left[  k\right]  \rightarrow\left[  n\right]  }%
\underbrace{\sum_{\sigma\in S_{n}}\left(  -1\right)  ^{\sigma}\prod_{i=1}%
^{k}a_{\kappa\left(  i\right)  ,\sigma\left(  \kappa\left(  i\right)  \right)
}}_{\substack{=0\\\text{(by (\ref{pf.prop.sol.powerdet.gen.a.1}))}}%
}=\sum_{\kappa:\left[  k\right]  \rightarrow\left[  n\right]  }0=0.
\end{align*}
This proves Proposition \ref{prop.sol.powerdet.gen.a}.
\end{proof}

\begin{remark}
Corollary \ref{cor.sol.powerdet.gen.prodrule}, Lemma
\ref{lem.sol.powerdet.gen.toofew} and Proposition
\ref{prop.sol.powerdet.gen.a} all remain valid if the commutative ring
$\mathbb{K}$ is replaced by a noncommutative ring $\mathbb{L}$, as long as we
use the conventions made in Section \ref{sect.sol.noncomm.polarization}. In
fact, the proofs given above still work when $\mathbb{K}$ is replaced by
$\mathbb{L}$, provided that we replace the reference to Lemma
\ref{lem.prodrule2} by a reference to Lemma \ref{lem.noncomm.prodrule2}.
\end{remark}

Part \textbf{(b)} of Exercise \ref{exe.powerdet.gen} is noticeably harder. We
prepare to it by studying permutations in $S_{n}$:

\begin{lemma}
\label{lem.sol.powerdet.gen.allbut1}Let $n\in\mathbb{N}$ and $p\in\left[
n\right]  $. Let $\alpha\in S_{n}$ and $\beta\in S_{n}$. Assume that%
\begin{equation}
\alpha\left(  i\right)  =\beta\left(  i\right)  \ \ \ \ \ \ \ \ \ \ \text{for
each }i\in\left[  n\right]  \setminus\left\{  p\right\}  .
\label{eq.lem.sol.powerdet.gen.allbut1.ass}%
\end{equation}
Then, $\alpha=\beta$.
\end{lemma}

\begin{vershort}
\begin{proof}
[Proof of Lemma \ref{lem.sol.powerdet.gen.allbut1}.]We have $\alpha\in S_{n}$.
In other words, $\alpha$ is a permutation of $\left[  n\right]  $ (since
$S_{n}$ is the set of all permutations of $\left\{  1,2,\ldots,n\right\}
=\left[  n\right]  $). In other words, $\alpha$ is a bijective map $\left[
n\right]  \rightarrow\left[  n\right]  $. The map $\alpha$ is bijective and
thus injective. Hence, Lemma \ref{lem.jectivity.pigeon0} \textbf{(c)} (applied
to $U=\left[  n\right]  $, $V=\left[  n\right]  $, $f=\alpha$ and $S=\left[
n\right]  \setminus\left\{  p\right\}  $) shows that
\[
\left\vert \alpha\left(  \left[  n\right]  \setminus\left\{  p\right\}
\right)  \right\vert =\left\vert \left[  n\right]  \setminus\left\{
p\right\}  \right\vert =\left\vert \left[  n\right]  \right\vert
-1\ \ \ \ \ \ \ \ \ \ \left(  \text{since }p\in\left[  n\right]  \right)  .
\]

Let $G$ denote the subset $\alpha\left(  \left[  n\right]  \setminus\left\{
p\right\}  \right)  $ of $\left[  n\right]  $. Then, $G=\alpha\left(  \left[
n\right]  \setminus\left\{  p\right\}  \right)  $, and thus $\left\vert
G\right\vert =\left\vert \alpha\left(  \left[  n\right]  \setminus\left\{
p\right\}  \right)  \right\vert =\left\vert \left[  n\right]  \right\vert -1$.

But $G$ is a subset of $\left[  n\right]  $; hence,%
\[
\left\vert \left[  n\right]  \setminus G\right\vert =\left\vert \left[
n\right]  \right\vert -\underbrace{\left\vert G\right\vert }_{=\left\vert
\left[  n\right]  \right\vert -1}=\left\vert \left[  n\right]  \right\vert
-\left(  \left\vert \left[  n\right]  \right\vert -1\right)  =1.
\]
In other words, $\left[  n\right]  \setminus G$ is a $1$-element set. Hence,
$\left[  n\right]  \setminus G=\left\{  q\right\}  $ for some element $q$.
Consider this $q$.

We have $q\in\left\{  q\right\}  =\left[  n\right]  \setminus G\subseteq
\left[  n\right]  $.

Now, $\alpha\left(  p\right)  \notin G$\ \ \ \ \footnote{\textit{Proof.}
Assume the contrary. Thus, $\alpha\left(  p\right)  \in G$. Hence,
$\alpha\left(  p\right)  \in G=\alpha\left(  \left[  n\right]  \setminus
\left\{  p\right\}  \right)  $. In other words, there exists an $i\in\left[
n\right]  \setminus\left\{  p\right\}  $ satisfying $\alpha\left(  p\right)
=\alpha\left(  i\right)  $. Consider this $i$.
\par
From $\alpha\left(  p\right)  =\alpha\left(  i\right)  $, we obtain $p=i$
(since $\alpha$ is injective). But $i\in\left[  n\right]  \setminus\left\{
p\right\}  $ shows that $i\notin\left\{  p\right\}  $, so that $i\neq p=i$.
This is clearly a contradiction. This contradiction proves that our assumption
was wrong. Qed.}. Combining this with $\alpha\left(  p\right)  \in\left[
n\right]  $, we obtain $\alpha\left(  p\right)  \in\left[  n\right]  \setminus
G=\left\{  q\right\}  $. Thus, $\alpha\left(  p\right)  =q$.

Also, $\beta\left(  p\right)  \notin G$\ \ \ \ \footnote{\textit{Proof.}
Assume the contrary. Thus, $\beta\left(  p\right)  \in G$. Hence,
$\beta\left(  p\right)  \in G=\alpha\left(  \left[  n\right]  \setminus
\left\{  p\right\}  \right)  $. In other words, there exists an $i\in\left[
n\right]  \setminus\left\{  p\right\}  $ satisfying $\beta\left(  p\right)
=\alpha\left(  i\right)  $. Consider this $i$.
\par
From (\ref{eq.lem.sol.powerdet.gen.allbut1.ass}), we obtain $\alpha\left(
i\right)  =\beta\left(  i\right)  $. Thus, $\beta\left(  p\right)
=\alpha\left(  i\right)  =\beta\left(  i\right)  $.
\par
But recall that the map $\alpha$ is injective. Similarly, the map $\beta$ is
injective. Thus, from $\beta\left(  p\right)  =\beta\left(  i\right)  $, we
obtain $p=i$. But $i\in\left[  n\right]  \setminus\left\{  p\right\}  $ shows
that $i\notin\left\{  p\right\}  $, so that $i\neq p=i$. This is a
contradiction. This contradiction proves that our assumption was wrong. Qed.}.
Combining this with $\beta\left(  p\right)  \in\left[  n\right]  $, we obtain
$\beta\left(  p\right)  \in\left[  n\right]  \setminus G=\left\{  q\right\}
$. Thus, $\beta\left(  p\right)  =q$. Comparing this with $\alpha\left(
p\right)  =q$, we obtain $\alpha\left(  p\right)  =\beta\left(  p\right)  $.

Now, each $i\in\left[  n\right]  $ satisfies $\alpha\left(  i\right)
=\beta\left(  i\right)  $\ \ \ \ \footnote{\textit{Proof.} Let $i\in\left[
n\right]  $. We must show that $\alpha\left(  i\right)  =\beta\left(
i\right)  $.
\par
If $i=p$, then this is clearly true (because $\alpha\left(  p\right)
=\beta\left(  p\right)  $). Hence, for the rest of this proof, we WLOG assume
that $i\neq p$. Hence, $i\in\left[  n\right]  \setminus\left\{  p\right\}  $.
Thus, (\ref{eq.lem.sol.powerdet.gen.allbut1.ass}) yields $\alpha\left(
i\right)  =\beta\left(  i\right)  $. Qed.}. In other words, $\alpha=\beta$.
This proves Lemma \ref{lem.sol.powerdet.gen.allbut1}.
\end{proof}
\end{vershort}

\begin{verlong}
\begin{proof}
[Proof of Lemma \ref{lem.sol.powerdet.gen.allbut1}.]We have $\alpha\in S_{n}$.
In other words, $\alpha$ is a permutation of $\left\{  1,2,\ldots,n\right\}  $
(since $S_{n}$ is the set of all permutations of $\left\{  1,2,\ldots
,n\right\}  $). In other words, $\alpha$ is a permutation of $\left[
n\right]  $ (since $\left\{  1,2,\ldots,n\right\}  =\left[  n\right]  $). In
other words, $\alpha$ is a bijective map $\left[  n\right]  \rightarrow\left[
n\right]  $. The map $\alpha$ is bijective and thus injective. Hence, Lemma
\ref{lem.jectivity.pigeon0} \textbf{(c)} (applied to $U=\left[  n\right]  $,
$V=\left[  n\right]  $, $f=\alpha$ and $S=\left[  n\right]  \setminus\left\{
p\right\}  $) shows that
\[
\left\vert \alpha\left(  \left[  n\right]  \setminus\left\{  p\right\}
\right)  \right\vert =\left\vert \left[  n\right]  \setminus\left\{
p\right\}  \right\vert =\left\vert \left[  n\right]  \right\vert
-1\ \ \ \ \ \ \ \ \ \ \left(  \text{since }p\in\left[  n\right]  \right)  .
\]

Let $G$ denote the subset $\alpha\left(  \left[  n\right]  \setminus\left\{
p\right\}  \right)  $ of $\left[  n\right]  $. Then, $G=\alpha\left(  \left[
n\right]  \setminus\left\{  p\right\}  \right)  $, and thus $\left\vert
G\right\vert =\left\vert \alpha\left(  \left[  n\right]  \setminus\left\{
p\right\}  \right)  \right\vert =\left\vert \left[  n\right]  \right\vert -1$.

But $G$ is a subset of $\left[  n\right]  $; hence,%
\[
\left\vert \left[  n\right]  \setminus G\right\vert =\left\vert \left[
n\right]  \right\vert -\underbrace{\left\vert G\right\vert }_{=\left\vert
\left[  n\right]  \right\vert -1}=\left\vert \left[  n\right]  \right\vert
-\left(  \left\vert \left[  n\right]  \right\vert -1\right)  =1.
\]
In other words, $\left[  n\right]  \setminus G$ is a $1$-element set. Hence,
$\left[  n\right]  \setminus G$ has the form $\left[  n\right]  \setminus
G=\left\{  q\right\}  $ for some element $q$. Consider this $q$.

We have $q\in\left\{  q\right\}  =\left[  n\right]  \setminus G\subseteq
\left[  n\right]  $.

Now, $\alpha\left(  p\right)  \notin G$\ \ \ \ \footnote{\textit{Proof.}
Assume the contrary. Thus, $\alpha\left(  p\right)  \in G$. Hence,
$\alpha\left(  p\right)  \in G=\alpha\left(  \left[  n\right]  \setminus
\left\{  p\right\}  \right)  $. In other words, there exists an $i\in\left[
n\right]  \setminus\left\{  p\right\}  $ satisfying $\alpha\left(  p\right)
=\alpha\left(  i\right)  $. Consider this $i$.
\par
From $i\in\left[  n\right]  \setminus\left\{  p\right\}  $, we obtain
$i\in\left[  n\right]  $ and $i\notin\left\{  p\right\}  $. From
$i\notin\left\{  p\right\}  $, we obtain $i\neq p$ and thus $p\neq i$.
\par
We have $p\in\left[  n\right]  $ and $i\in\left[  n\right]  \setminus\left\{
p\right\}  \subseteq\left[  n\right]  $. But the map $\alpha$ is injective. In
other words, if $u\in\left[  n\right]  $ and $v\in\left[  n\right]  $ satisfy
$\alpha\left(  u\right)  =\alpha\left(  v\right)  $, then $u=v$. Applying this
to $u=p$ and $v=i$, we obtain $p=i$ (since $\alpha\left(  p\right)
=\alpha\left(  i\right)  $). This contradicts $p\neq i$. This contradiction
proves that our assumption was wrong. Qed.}. Combining this with
$\alpha\left(  p\right)  \in\left[  n\right]  $, we obtain $\alpha\left(
p\right)  \in\left[  n\right]  \setminus G=\left\{  q\right\}  $. Thus,
$\alpha\left(  p\right)  =q$.

Also, $\beta\left(  p\right)  \notin G$\ \ \ \ \footnote{\textit{Proof.}
Assume the contrary. Thus, $\beta\left(  p\right)  \in G$. Hence,
$\beta\left(  p\right)  \in G=\alpha\left(  \left[  n\right]  \setminus
\left\{  p\right\}  \right)  $. In other words, there exists an $i\in\left[
n\right]  \setminus\left\{  p\right\}  $ satisfying $\beta\left(  p\right)
=\alpha\left(  i\right)  $. Consider this $i$.
\par
From $i\in\left[  n\right]  \setminus\left\{  p\right\}  $, we obtain
$i\in\left[  n\right]  $ and $i\notin\left\{  p\right\}  $. From
$i\notin\left\{  p\right\}  $, we obtain $i\neq p$ and thus $p\neq i$.
\par
From (\ref{eq.lem.sol.powerdet.gen.allbut1.ass}), we obtain $\alpha\left(
i\right)  =\beta\left(  i\right)  $. Thus, $\beta\left(  p\right)
=\alpha\left(  i\right)  =\beta\left(  i\right)  $.
\par
But recall that $\alpha$ is an injective map $\left[  n\right]  \rightarrow
\left[  n\right]  $. The same argument (applied to $\beta$ instead of $\alpha
$) shows that $\beta$ is an injective map $\left[  n\right]  \rightarrow
\left[  n\right]  $. In particular, $\beta$ is injective. In other words, if
$u\in\left[  n\right]  $ and $v\in\left[  n\right]  $ satisfy $\beta\left(
u\right)  =\beta\left(  v\right)  $, then $u=v$. Applying this to $u=p$ and
$v=i$, we obtain $p=i$ (since $\beta\left(  p\right)  =\beta\left(  i\right)
$). This contradicts $p\neq i$. This contradiction proves that our assumption
was wrong. Qed.}. Combining this with $\beta\left(  p\right)  \in\left[
n\right]  $, we obtain $\beta\left(  p\right)  \in\left[  n\right]  \setminus
G=\left\{  q\right\}  $. Thus, $\beta\left(  p\right)  =q$. Comparing this
with $\alpha\left(  p\right)  =q$, we obtain $\alpha\left(  p\right)
=\beta\left(  p\right)  $.

Now, each $i\in\left[  n\right]  $ satisfies $\alpha\left(  i\right)
=\beta\left(  i\right)  $\ \ \ \ \footnote{\textit{Proof.} Let $i\in\left[
n\right]  $. We must show that $\alpha\left(  i\right)  =\beta\left(
i\right)  $.
\par
If $i=p$, then this is clearly true (because $\alpha\left(  p\right)
=\beta\left(  p\right)  $). Hence, for the rest of this proof, we can WLOG
assume that we don't have $i=p$. Assume this.
\par
We have $i\neq p$ (since we don't have $i=p$). Combining this with
$i\in\left[  n\right]  $, we obtain $i\in\left[  n\right]  \setminus\left\{
p\right\}  $. Thus, (\ref{eq.lem.sol.powerdet.gen.allbut1.ass}) yields
$\alpha\left(  i\right)  =\beta\left(  i\right)  $. Qed.}. In other words,
$\alpha=\beta$. This proves Lemma \ref{lem.sol.powerdet.gen.allbut1}.
\end{proof}
\end{verlong}

\begin{lemma}
\label{lem.sol.powerdet.gen.kappas}Let $n\geq1$ be an integer. Then, the map%
\begin{align*}
S_{n}  &  \rightarrow\left\{  f:\left[  n-1\right]  \rightarrow\left[
n\right]  \ \mid\ \left\vert f\left(  \left[  n-1\right]  \right)  \right\vert
\geq n-1\right\}  ,\\
\tau &  \mapsto\tau\mid_{\left[  n-1\right]  }%
\end{align*}
is well-defined and bijective.
\end{lemma}

\begin{proof}
[Proof of Lemma \ref{lem.sol.powerdet.gen.kappas}.]We have $n-1\in\mathbb{N}$
(since $n\geq1$) and thus $\left\vert \left[  n-1\right]  \right\vert =n-1$.

\begin{vershort}
Also, $n\in\left[  n\right]  $ (since $n\geq1$) and $\left[  n-1\right]
=\left[  n\right]  \setminus\left\{  n\right\}  \subseteq\left[  n\right]  $.
\end{vershort}

\begin{verlong}
Also, $n\in\left[  n\right]  $ (since $n\geq1$) and%
\[
\underbrace{\left[  n\right]  }_{=\left\{  1,2,\ldots,n\right\}  }%
\setminus\left\{  n\right\}  =\left\{  1,2,\ldots,n\right\}  \setminus\left\{
n\right\}  =\left\{  1,2,\ldots,n-1\right\}  .
\]
Thus, $\left[  n-1\right]  =\left\{  1,2,\ldots,n-1\right\}  =\left[
n\right]  \setminus\left\{  n\right\}  \subseteq\left[  n\right]  $. In other
words, $\left[  n-1\right]  $ is a subset of $\left[  n\right]  $.
\end{verlong}

Define a set $Y$ by%
\begin{equation}
Y=\left\{  f:\left[  n-1\right]  \rightarrow\left[  n\right]  \ \mid
\ \left\vert f\left(  \left[  n-1\right]  \right)  \right\vert \geq
n-1\right\}  . \label{pf.lem.sol.powerdet.gen.kappas.Y=}%
\end{equation}

\begin{vershort}
Each $\tau\in S_{n}$ satisfies $\tau\mid_{\left[  n-1\right]  }\in
Y$\ \ \ \ \footnote{\textit{Proof.} Let $\tau\in S_{n}$.
\par
We have $\tau\in S_{n}$. In other words, $\tau$ is a permutation of $\left[
n\right]  $ (since $S_{n}$ is the set of all permutations of $\left\{
1,2,\ldots,n\right\}  =\left[  n\right]  $). In other words, $\tau$ is a
bijective map $\left[  n\right]  \rightarrow\left[  n\right]  $. The map
$\tau$ is bijective and thus injective. Hence, Lemma
\ref{lem.jectivity.pigeon0} \textbf{(c)} (applied to $U=\left[  n\right]  $,
$V=\left[  n\right]  $, $f=\tau$ and $S=\left[  n-1\right]  $) shows that
$\left\vert \tau\left(  \left[  n-1\right]  \right)  \right\vert =\left\vert
\left[  n-1\right]  \right\vert =n-1$.
\par
Now, $\tau\mid_{\left[  n-1\right]  }$ is a map $\left[  n-1\right]
\rightarrow\left[  n\right]  $ and satisfies $\left\vert \underbrace{\left(
\tau\mid_{\left[  n-1\right]  }\right)  \left(  \left[  n-1\right]  \right)
}_{=\tau\left(  \left[  n-1\right]  \right)  }\right\vert =\left\vert
\tau\left(  \left[  n-1\right]  \right)  \right\vert =n-1\geq n-1$. Hence,
$\tau\mid_{\left[  n-1\right]  }$ is a map $f:\left[  n-1\right]
\rightarrow\left[  n\right]  $ satisfying $\left\vert f\left(  \left[
n-1\right]  \right)  \right\vert \geq n-1$. In other words,%
\[
\tau\mid_{\left[  n-1\right]  }\in\left\{  f:\left[  n-1\right]
\rightarrow\left[  n\right]  \ \mid\ \left\vert f\left(  \left[  n-1\right]
\right)  \right\vert \geq n-1\right\}  .
\]
In light of (\ref{pf.lem.sol.powerdet.gen.kappas.Y=}), this rewrites as
$\tau\mid_{\left[  n-1\right]  }\in Y$. Qed.}. Hence, we can define a map
$T:S_{n}\rightarrow Y$ by%
\[
\left(  T\left(  \tau\right)  =\tau\mid_{\left[  n-1\right]  }%
\ \ \ \ \ \ \ \ \ \ \text{for all }\tau\in S_{n}\right)  .
\]
Consider this map $T$.
\end{vershort}

\begin{verlong}
Each $\tau\in S_{n}$ satisfies $\tau\mid_{\left[  n-1\right]  }\in
Y$\ \ \ \ \footnote{\textit{Proof.} Let $\tau\in S_{n}$.
\par
We have $\tau\in S_{n}$. In other words, $\tau$ is a permutation of $\left\{
1,2,\ldots,n\right\}  $ (since $S_{n}$ is the set of all permutations of
$\left\{  1,2,\ldots,n\right\}  $). In other words, $\tau$ is a permutation of
$\left[  n\right]  $ (since $\left\{  1,2,\ldots,n\right\}  =\left[  n\right]
$). In other words, $\tau$ is a bijective map $\left[  n\right]
\rightarrow\left[  n\right]  $. The map $\tau$ is bijective and thus
injective. Hence, Lemma \ref{lem.jectivity.pigeon0} \textbf{(c)} (applied to
$U=\left[  n\right]  $, $V=\left[  n\right]  $, $f=\tau$ and $S=\left[
n-1\right]  $) shows that $\left\vert \tau\left(  \left[  n-1\right]  \right)
\right\vert =\left\vert \left[  n-1\right]  \right\vert =n-1$.
\par
Now, $\tau\mid_{\left[  n-1\right]  }$ is a map $\left[  n-1\right]
\rightarrow\left[  n\right]  $ and satisfies $\left\vert \underbrace{\left(
\tau\mid_{\left[  n-1\right]  }\right)  \left(  \left[  n-1\right]  \right)
}_{=\tau\left(  \left[  n-1\right]  \right)  }\right\vert =\left\vert
\tau\left(  \left[  n-1\right]  \right)  \right\vert =n-1\geq n-1$. Hence,
$\tau\mid_{\left[  n-1\right]  }$ is a map $f:\left[  n-1\right]
\rightarrow\left[  n\right]  $ satisfying $\left\vert f\left(  \left[
n-1\right]  \right)  \right\vert \geq n-1$. In other words,%
\[
\tau\mid_{\left[  n-1\right]  }\in\left\{  f:\left[  n-1\right]
\rightarrow\left[  n\right]  \ \mid\ \left\vert f\left(  \left[  n-1\right]
\right)  \right\vert \geq n-1\right\}  .
\]
In light of (\ref{pf.lem.sol.powerdet.gen.kappas.Y=}), this rewrites as
$\tau\mid_{\left[  n-1\right]  }\in Y$. Qed.}. Hence, we can define a map
$T:S_{n}\rightarrow Y$ by%
\[
\left(  T\left(  \tau\right)  =\tau\mid_{\left[  n-1\right]  }%
\ \ \ \ \ \ \ \ \ \ \text{for all }\tau\in S_{n}\right)  .
\]
Consider this map $T$.
\end{verlong}

\begin{vershort}
Each $g\in Y$ satisfies $g\in T\left(  S_{n}\right)  $%
\ \ \ \ \footnote{\textit{Proof.} Let $g\in Y$. Thus, $g\in Y=\left\{
f:\left[  n-1\right]  \rightarrow\left[  n\right]  \ \mid\ \left\vert f\left(
\left[  n-1\right]  \right)  \right\vert \geq n-1\right\}  $. In other words,
$g$ is a map $\left[  n-1\right]  \rightarrow\left[  n\right]  $ and satisfies
$\left\vert g\left(  \left[  n-1\right]  \right)  \right\vert \geq n-1$.
\par
Thus, $\left\vert g\left(  \left[  n-1\right]  \right)  \right\vert \geq
n-1=\left\vert \left[  n-1\right]  \right\vert $. Hence, Lemma
\ref{lem.jectivity.pigeon0} \textbf{(b)} (applied to $U=\left[  n-1\right]  $,
$V=\left[  n\right]  $ and $f=g$) shows that the map $g$ is injective. Hence,
Lemma \ref{lem.jectivity.pigeon0} \textbf{(c)} (applied to $U=\left[
n-1\right]  $, $V=\left[  n\right]  $, $f=g$ and $S=\left[  n-1\right]  $)
shows that $\left\vert g\left(  \left[  n-1\right]  \right)  \right\vert
=\left\vert \left[  n-1\right]  \right\vert =n-1$. Since $g\left(  \left[
n-1\right]  \right)  $ is a subset of $\left[  n\right]  $, we have%
\[
\left\vert \left[  n\right]  \setminus g\left(  \left[  n-1\right]  \right)
\right\vert =\underbrace{\left\vert \left[  n\right]  \right\vert }%
_{=n}-\underbrace{\left\vert g\left(  \left[  n-1\right]  \right)  \right\vert
}_{=n-1}=n-\left(  n-1\right)  =1.
\]
In other words, $\left[  n\right]  \setminus g\left(  \left[  n-1\right]
\right)  $ is a $1$-element set. In other words, $\left[  n\right]  \setminus
g\left(  \left[  n-1\right]  \right)  =\left\{  p\right\}  $ for some element
$p$. Consider this $p$.
\par
We have $p\in\left\{  p\right\}  =\left[  n\right]  \setminus g\left(  \left[
n-1\right]  \right)  \subseteq\left[  n\right]  $.
\par
Now, define a map $\sigma:\left[  n\right]  \rightarrow\left[  n\right]  $ by%
\[
\left(  \sigma\left(  i\right)  =%
\begin{cases}
g\left(  i\right)  , & \text{if }i\in\left[  n-1\right]  ;\\
p, & \text{if }i\notin\left[  n-1\right]
\end{cases}
\ \ \ \ \ \ \ \ \ \ \text{for each }i\in\left[  n\right]  \right)  .
\]
(This is well-defined, because each $i\in\left[  n\right]  $ satisfies $%
\begin{cases}
g\left(  i\right)  , & \text{if }i\in\left[  n-1\right]  ;\\
p, & \text{if }i\notin\left[  n-1\right]
\end{cases}
\in\left[  n\right]  $.)
\par
We shall now show that the map $\sigma$ is surjective.
\par
Indeed, $n\notin\left\{  1,2,\ldots,n-1\right\}  =\left[  n-1\right]  $ but
$n\in\left[  n\right]  $. Hence, the definition of $\sigma$ yields%
\[
\sigma\left(  n\right)  =%
\begin{cases}
g\left(  n\right)  , & \text{if }n\in\left[  n-1\right]  ;\\
p, & \text{if }n\notin\left[  n-1\right]
\end{cases}
=p\ \ \ \ \ \ \ \ \ \ \left(  \text{since }n\notin\left[  n-1\right]  \right)
.
\]
But $\sigma\left(  \underbrace{n}_{\in\left[  n\right]  }\right)  \in
\sigma\left(  \left[  n\right]  \right)  $, so that $p=\sigma\left(  n\right)
\in\sigma\left(  \left[  n\right]  \right)  $.
\par
Also, each $i\in\left[  n-1\right]  $ satisfies%
\begin{equation}
\sigma\left(  i\right)  =%
\begin{cases}
g\left(  i\right)  , & \text{if }i\in\left[  n-1\right]  ;\\
p, & \text{if }i\notin\left[  n-1\right]
\end{cases}
=g\left(  i\right)  \ \ \ \ \ \ \ \ \ \ \left(  \text{since }i\in\left[
n-1\right]  \right)  . \label{pf.lem.sol.powerdet.gen.kappas.short.sur.24}%
\end{equation}
\par
Now,%
\[
g\left(  \left[  n-1\right]  \right)  =\left\{  \underbrace{g\left(  i\right)
}_{\substack{=\sigma\left(  i\right)  \\\text{(by
(\ref{pf.lem.sol.powerdet.gen.kappas.short.sur.24}))}}}\ \mid\ i\in\left[
n-1\right]  \right\}  =\left\{  \sigma\left(  i\right)  \ \mid\ i\in\left[
n-1\right]  \right\}  =\sigma\left(  \underbrace{\left[  n-1\right]
}_{\subseteq\left[  n\right]  }\right)  \subseteq\sigma\left(  \left[
n\right]  \right)  .
\]
\par
If $X$ is a set, and if $Y$ is a subset of $X$, then $X=Y\cup\left(
X\setminus Y\right)  $. Applying this to $X=\left[  n\right]  $ and
$Y=g\left(  \left[  n-1\right]  \right)  $, we obtain%
\[
\left[  n\right]  =\underbrace{g\left(  \left[  n-1\right]  \right)
}_{\subseteq\sigma\left(  \left[  n\right]  \right)  }\cup\underbrace{\left(
\left[  n\right]  \setminus g\left(  \left[  n-1\right]  \right)  \right)
}_{\substack{=\left\{  p\right\}  \subseteq\sigma\left(  \left[  n\right]
\right)  \\\text{(since }p\in\sigma\left(  \left[  n\right]  \right)
\text{)}}}\subseteq\sigma\left(  \left[  n\right]  \right)  \cup\sigma\left(
\left[  n\right]  \right)  =\sigma\left(  \left[  n\right]  \right)  .
\]
Combining this with the obvious relation $\sigma\left(  \left[  n\right]
\right)  \subseteq\left[  n\right]  $, we obtain $\left[  n\right]
=\sigma\left(  \left[  n\right]  \right)  $. In other words, the map $\sigma$
is surjective.
\par
Lemma \ref{lem.jectivity.pigeon-surj} (applied to $U=\left[  n\right]  $,
$V=\left[  n\right]  $ and $f=\sigma$) shows that we have the following
logical equivalence:%
\[
\left(  \sigma\text{ is surjective}\right)  \ \Longleftrightarrow\ \left(
\sigma\text{ is bijective}\right)
\]
(since $\left\vert \left[  n\right]  \right\vert \leq\left\vert \left[
n\right]  \right\vert $). Thus, $\sigma$ is bijective (since $\sigma$ is
surjective). Hence, $\sigma$ is a bijective map $\left[  n\right]
\rightarrow\left[  n\right]  $. In other words, $\sigma$ is a permutation of
$\left[  n\right]  $. In other words, $\sigma\in S_{n}$ (since $S_{n}$ is the
set of all permutations of $\left\{  1,2,\ldots,n\right\}  =\left[  n\right]
$). Hence, $T\left(  \sigma\right)  $ is well-defined.
\par
The definition of $T$ yields $T\left(  \sigma\right)  =\sigma\mid_{\left[
n-1\right]  }$. Thus, each $i\in\left[  n-1\right]  $ satisfies%
\[
\left(  T\left(  \sigma\right)  \right)  \left(  i\right)  =\left(  \sigma
\mid_{\left[  n-1\right]  }\right)  \left(  i\right)  =\sigma\left(  i\right)
=g\left(  i\right)  \ \ \ \ \ \ \ \ \ \ \left(  \text{by
(\ref{pf.lem.sol.powerdet.gen.kappas.short.sur.24})}\right)  .
\]
In other words, $T\left(  \sigma\right)  =g$. Thus, $g=T\left(
\underbrace{\sigma}_{\in S_{n}}\right)  \in T\left(  S_{n}\right)  $. Qed.}.
In other words, $Y\subseteq T\left(  S_{n}\right)  $. Combining this with
$T\left(  S_{n}\right)  \subseteq Y$ (which is obvious), we obtain $Y=T\left(
S_{n}\right)  $. In other words, the map $T$ is surjective.
\end{vershort}

\begin{verlong}
Each $g\in Y$ satisfies $g\in T\left(  S_{n}\right)  $%
\ \ \ \ \footnote{\textit{Proof.} Let $g\in Y$. Thus, $g\in Y=\left\{
f:\left[  n-1\right]  \rightarrow\left[  n\right]  \ \mid\ \left\vert f\left(
\left[  n-1\right]  \right)  \right\vert \geq n-1\right\}  $. In other words,
$g$ is a map $f:\left[  n-1\right]  \rightarrow\left[  n\right]  $ satisfying
$\left\vert f\left(  \left[  n-1\right]  \right)  \right\vert \geq n-1$. In
other words, $g$ is a map $\left[  n-1\right]  \rightarrow\left[  n\right]  $
and satisfies $\left\vert g\left(  \left[  n-1\right]  \right)  \right\vert
\geq n-1$.
\par
Thus, $\left\vert g\left(  \left[  n-1\right]  \right)  \right\vert \geq
n-1=\left\vert \left[  n-1\right]  \right\vert $. Hence, Lemma
\ref{lem.jectivity.pigeon0} \textbf{(b)} (applied to $U=\left[  n-1\right]  $,
$V=\left[  n\right]  $ and $f=g$) shows that the map $g$ is injective. Hence,
Lemma \ref{lem.jectivity.pigeon0} \textbf{(c)} (applied to $U=\left[
n-1\right]  $, $V=\left[  n\right]  $, $f=g$ and $S=\left[  n-1\right]  $)
shows that $\left\vert g\left(  \left[  n-1\right]  \right)  \right\vert
=\left\vert \left[  n-1\right]  \right\vert =n-1$. Since $g\left(  \left[
n-1\right]  \right)  $ is a subset of $\left[  n\right]  $, we have%
\[
\left\vert \left[  n\right]  \setminus g\left(  \left[  n-1\right]  \right)
\right\vert =\underbrace{\left\vert \left[  n\right]  \right\vert }%
_{=n}-\underbrace{\left\vert g\left(  \left[  n-1\right]  \right)  \right\vert
}_{=n-1}=n-\left(  n-1\right)  =1.
\]
In other words, $\left[  n\right]  \setminus g\left(  \left[  n-1\right]
\right)  $ is a $1$-element set. In other words, $\left[  n\right]  \setminus
g\left(  \left[  n-1\right]  \right)  =\left\{  p\right\}  $ for some element
$p$. Consider this $p$.
\par
We have $p\in\left\{  p\right\}  =\left[  n\right]  \setminus g\left(  \left[
n-1\right]  \right)  \subseteq\left[  n\right]  $.
\par
Now, each $i\in\left[  n\right]  $ satisfies $%
\begin{cases}
g\left(  i\right)  , & \text{if }i\in\left[  n-1\right]  ;\\
p, & \text{if }i\notin\left[  n-1\right]
\end{cases}
\in\left[  n\right]  $ (because we have $g\left(  i\right)  \in\left[
n\right]  $ whenever $i\in\left[  n-1\right]  $ (since $g$ is a map $\left[
n-1\right]  \rightarrow\left[  n\right]  $), and because we have $p\in\left[
n\right]  $ whenever $i\notin\left[  n-1\right]  $ (since $p\in\left[
n\right]  $)).
\par
Now, define a map $\sigma:\left[  n\right]  \rightarrow\left[  n\right]  $ by%
\[
\left(  \sigma\left(  i\right)  =%
\begin{cases}
g\left(  i\right)  , & \text{if }i\in\left[  n-1\right]  ;\\
p, & \text{if }i\notin\left[  n-1\right]
\end{cases}
\ \ \ \ \ \ \ \ \ \ \text{for each }i\in\left[  n\right]  \right)  .
\]
(This is well-defined, because each $i\in\left[  n\right]  $ satisfies $%
\begin{cases}
g\left(  i\right)  , & \text{if }i\in\left[  n-1\right]  ;\\
p, & \text{if }i\notin\left[  n-1\right]
\end{cases}
\in\left[  n\right]  $.)
\par
We shall now show that the map $\sigma$ is surjective.
\par
Indeed, $n\notin\left\{  1,2,\ldots,n-1\right\}  =\left[  n-1\right]  $ but
$n\in\left[  n\right]  $. Hence, the definition of $\sigma$ yields%
\[
\sigma\left(  n\right)  =%
\begin{cases}
g\left(  n\right)  , & \text{if }n\in\left[  n-1\right]  ;\\
p, & \text{if }n\notin\left[  n-1\right]
\end{cases}
=p\ \ \ \ \ \ \ \ \ \ \left(  \text{since }n\notin\left[  n-1\right]  \right)
.
\]
But $\sigma\left(  \underbrace{n}_{\in\left[  n\right]  }\right)  \in
\sigma\left(  \left[  n\right]  \right)  $, so that $p=\sigma\left(  n\right)
\in\sigma\left(  \left[  n\right]  \right)  $.
\par
Also, each $i\in\left[  n-1\right]  $ satisfies $i\in\left[  n-1\right]
\subseteq\left[  n\right]  $. Therefore, for each $i\in\left[  n-1\right]  $,
the value $\sigma\left(  i\right)  $ is well-defined. Furthermore, each
$i\in\left[  n-1\right]  $ satisfies%
\begin{equation}
\sigma\left(  i\right)  =%
\begin{cases}
g\left(  i\right)  , & \text{if }i\in\left[  n-1\right]  ;\\
p, & \text{if }i\notin\left[  n-1\right]
\end{cases}
=g\left(  i\right)  \ \ \ \ \ \ \ \ \ \ \left(  \text{since }i\in\left[
n-1\right]  \right)  . \label{pf.lem.sol.powerdet.gen.kappas.sur.24}%
\end{equation}
\par
Now,%
\[
g\left(  \left[  n-1\right]  \right)  =\left\{  \underbrace{g\left(  i\right)
}_{\substack{=\sigma\left(  i\right)  \\\text{(by
(\ref{pf.lem.sol.powerdet.gen.kappas.sur.24}))}}}\ \mid\ i\in\left[
n-1\right]  \right\}  =\left\{  \sigma\left(  i\right)  \ \mid\ i\in\left[
n-1\right]  \right\}  =\sigma\left(  \underbrace{\left[  n-1\right]
}_{\subseteq\left[  n\right]  }\right)  \subseteq\sigma\left(  \left[
n\right]  \right)  .
\]
\par
If $X$ is a set, and if $Y$ is a subset of $X$, then $X=Y\cup\left(
X\setminus Y\right)  $. Applying this to $X=\left[  n\right]  $ and
$Y=g\left(  \left[  n-1\right]  \right)  $, we obtain%
\[
\left[  n\right]  =\underbrace{g\left(  \left[  n-1\right]  \right)
}_{\subseteq\sigma\left(  \left[  n\right]  \right)  }\cup\underbrace{\left(
\left[  n\right]  \setminus g\left(  \left[  n-1\right]  \right)  \right)
}_{\substack{=\left\{  p\right\}  \subseteq\sigma\left(  \left[  n\right]
\right)  \\\text{(since }p\in\sigma\left(  \left[  n\right]  \right)
\text{)}}}\subseteq\sigma\left(  \left[  n\right]  \right)  \cup\sigma\left(
\left[  n\right]  \right)  =\sigma\left(  \left[  n\right]  \right)  .
\]
Combining this with the obvious relation $\sigma\left(  \left[  n\right]
\right)  \subseteq\left[  n\right]  $, we obtain $\left[  n\right]
=\sigma\left(  \left[  n\right]  \right)  $. In other words, the map $\sigma$
is surjective.
\par
Lemma \ref{lem.jectivity.pigeon-surj} (applied to $U=\left[  n\right]  $,
$V=\left[  n\right]  $ and $f=\sigma$) shows that we have the following
logical equivalence:%
\[
\left(  \sigma\text{ is surjective}\right)  \ \Longleftrightarrow\ \left(
\sigma\text{ is bijective}\right)
\]
(since $\left\vert \left[  n\right]  \right\vert \leq\left\vert \left[
n\right]  \right\vert $). Thus, $\sigma$ is bijective (since $\sigma$ is
surjective). Hence, $\sigma$ is a bijective map $\left[  n\right]
\rightarrow\left[  n\right]  $. In other words, $\sigma$ is a permutation of
$\left[  n\right]  $. In other words, $\sigma$ is a permutation of $\left\{
1,2,\ldots,n\right\}  $ (since $\left[  n\right]  =\left\{  1,2,\ldots
,n\right\}  $). In other words, $\sigma\in S_{n}$ (since $S_{n}$ is the set of
all permutations of $\left\{  1,2,\ldots,n\right\}  $). Hence, $T\left(
\sigma\right)  $ is well-defined.
\par
The definition of $T$ yields $T\left(  \sigma\right)  =\sigma\mid_{\left[
n-1\right]  }$. Thus, each $i\in\left[  n-1\right]  $ satisfies%
\[
\underbrace{\left(  T\left(  \sigma\right)  \right)  }_{=\sigma\mid_{\left[
n-1\right]  }}\left(  i\right)  =\left(  \sigma\mid_{\left[  n-1\right]
}\right)  \left(  i\right)  =\sigma\left(  i\right)  =g\left(  i\right)
\ \ \ \ \ \ \ \ \ \ \left(  \text{by
(\ref{pf.lem.sol.powerdet.gen.kappas.sur.24})}\right)  .
\]
In other words, $T\left(  \sigma\right)  =g$. Thus, $g=T\left(
\underbrace{\sigma}_{\in S_{n}}\right)  \in T\left(  S_{n}\right)  $. Qed.}.
In other words, $Y\subseteq T\left(  S_{n}\right)  $. Combining this with
$T\left(  S_{n}\right)  \subseteq Y$ (which is obvious), we obtain $Y=T\left(
S_{n}\right)  $. In other words, the map $T$ is surjective.
\end{verlong}

If $\alpha\in S_{n}$ and $\beta\in S_{n}$ satisfy $T\left(  \alpha\right)
=T\left(  \beta\right)  $, then $\alpha=\beta$%
\ \ \ \ \footnote{\textit{Proof.} Let $\alpha\in S_{n}$ and $\beta\in S_{n}$
be such that $T\left(  \alpha\right)  =T\left(  \beta\right)  $. We must show
that $\alpha=\beta$.
\par
The definition of $T$ yields $T\left(  \alpha\right)  =\alpha\mid_{\left[
n-1\right]  }$ and $T\left(  \beta\right)  =\beta\mid_{\left[  n-1\right]  }$.
Hence, $\alpha\mid_{\left[  n-1\right]  }=T\left(  \alpha\right)  =T\left(
\beta\right)  =\beta\mid_{\left[  n-1\right]  }$.
\par
Now, each $i\in\left[  n-1\right]  $ satisfies $\alpha\left(  i\right)
=\underbrace{\left(  \alpha\mid_{\left[  n-1\right]  }\right)  }_{=\beta
\mid_{\left[  n-1\right]  }}\left(  i\right)  =\left(  \beta\mid_{\left[
n-1\right]  }\right)  \left(  i\right)  =\beta\left(  i\right)  $. Thus,
$\alpha\left(  i\right)  =\beta\left(  i\right)  $ for each $i\in\left[
n-1\right]  $. In other words, $\alpha\left(  i\right)  =\beta\left(
i\right)  $ for each $i\in\left[  n\right]  \setminus\left\{  n\right\}  $
(since $\left[  n-1\right]  =\left[  n\right]  \setminus\left\{  n\right\}
$). Thus, Lemma \ref{lem.sol.powerdet.gen.allbut1} (applied to $p=n$) yields
$\alpha=\beta$. Qed.}. In other words, the map $T$ is injective.

The map $T$ is surjective and injective. In other words, the map $T$ is bijective.

Recall that $T$ is a map from $S_{n}$ to $Y$. In other words, $T$ is a map
from $S_{n}$ to $\left\{  f:\left[  n-1\right]  \rightarrow\left[  n\right]
\ \mid\ \left\vert f\left(  \left[  n-1\right]  \right)  \right\vert \geq
n-1\right\}  $ (since \newline$Y=\left\{  f:\left[  n-1\right]  \rightarrow
\left[  n\right]  \ \mid\ \left\vert f\left(  \left[  n-1\right]  \right)
\right\vert \geq n-1\right\}  $). Hence, the map
\begin{align*}
S_{n}  &  \rightarrow\left\{  f:\left[  n-1\right]  \rightarrow\left[
n\right]  \ \mid\ \left\vert f\left(  \left[  n-1\right]  \right)  \right\vert
\geq n-1\right\}  ,\\
\tau &  \mapsto\tau\mid_{\left[  n-1\right]  }%
\end{align*}
is precisely the map $T$ (since $T\left(  \tau\right)  =\tau\mid_{\left[
n-1\right]  }$ for all $\tau\in S_{n}$). Hence, this map is well-defined and
bijective (since $T$ is bijective). This proves Lemma
\ref{lem.sol.powerdet.gen.kappas}.
\end{proof}

\begin{lemma}
\label{lem.sol.powerdet.gen.taures}Let $n\geq1$ be an integer. Let $\tau\in
S_{n}$. Then,%
\[
\sum_{\sigma\in S_{n}}\left(  -1\right)  ^{\sigma}\prod_{i=1}^{n-1}a_{\left(
\tau\mid_{\left[  n-1\right]  }\right)  \left(  i\right)  ,\sigma\left(
\left(  \tau\mid_{\left[  n-1\right]  }\right)  \left(  i\right)  \right)
}=\sum_{q=1}^{n}\left(  -1\right)  ^{\tau\left(  n\right)  +q}\det\left(
A_{\sim\left(  \tau\left(  n\right)  \right)  ,\sim q}\right)  .
\]

\end{lemma}

\begin{vershort}
\begin{proof}
[Proof of Lemma \ref{lem.sol.powerdet.gen.taures}.]We have $\tau\in S_{n}$. In
other words, $\tau$ is a permutation of $\left[  n\right]  $ (since $S_{n}$ is
the set of all permutations of $\left\{  1,2,\ldots,n\right\}  =\left[
n\right]  $). In other words, $\tau$ is a bijective map $\left[  n\right]
\rightarrow\left[  n\right]  $. Thus, $\tau:\left[  n\right]  \rightarrow
\left[  n\right]  $ is a bijection.

Also, $n\in\left[  n\right]  $ (since $n\geq1$) and $\left[  n-1\right]
=\left[  n\right]  \setminus\left\{  n\right\}  \subseteq\left[  n\right]  $.

Let $p=\tau\left(  n\right)  $. Then, $p=\tau\left(  n\right)  \in\left[
n\right]  =\left\{  1,2,\ldots,n\right\}  $. From $p=\tau\left(  n\right)  $,
we obtain $n=\tau^{-1}\left(  p\right)  $.

If $b_{1},b_{2},\ldots,b_{n}$ are $n$ elements of $\mathbb{K}$, then%
\begin{equation}
\prod_{i=1}^{n-1}b_{\left(  \tau\mid_{\left[  n-1\right]  }\right)  \left(
i\right)  }=\prod_{\substack{i\in\left\{  1,2,\ldots,n\right\}  ;\\i\neq
p}}b_{i} \label{pf.lem.sol.powerdet.gen.taures.short.bprod}%
\end{equation}
\footnote{\textit{Proof of (\ref{pf.lem.sol.powerdet.gen.taures.short.bprod}%
):} Let $b_{1},b_{2},\ldots,b_{n}$ be $n$ elements of $\mathbb{K}$.
\par
For each $i\in\left[  n\right]  $, we have the following chain of logical
equivalences:%
\[
\left(  i\neq\underbrace{n}_{=\tau^{-1}\left(  p\right)  }\right)
\ \Longleftrightarrow\ \left(  i\neq\tau^{-1}\left(  p\right)  \right)
\ \Longleftrightarrow\ \left(  \tau\left(  i\right)  \neq p\right)
\]
(since $\tau$ is a bijection). In other words, for each $i\in\left[  n\right]
$, the condition $\left(  i\neq n\right)  $ holds if and only if the condition
$\left(  \tau\left(  i\right)  \neq p\right)  $ holds. Thus, we have the
following equality of product signs:%
\[
\prod_{\substack{i\in\left[  n\right]  ;\\i\neq n}}=\prod_{\substack{i\in
\left[  n\right]  ;\\\tau\left(  i\right)  \neq p}}.
\]
Hence,%
\begin{align}
\prod_{\substack{i\in\left[  n\right]  ;\\\tau\left(  i\right)  \neq p}}  &
=\prod_{\substack{i\in\left[  n\right]  ;\\i\neq n}}=\prod_{i\in\left[
n\right]  \setminus\left\{  n\right\}  }=\prod_{i\in\left[  n-1\right]
}\ \ \ \ \ \ \ \ \ \ \left(  \text{since }\left[  n\right]  \setminus\left\{
n\right\}  =\left[  n-1\right]  \right) \nonumber\\
&  =\prod_{i=1}^{n-1}. \label{pf.lem.sol.powerdet.gen.taures.short.bprod.preq}%
\end{align}
\par
Now, the map $\tau:\left[  n\right]  \rightarrow\left[  n\right]  $ is a
bijection. Thus, we can substitute $\tau\left(  i\right)  $ for $i$ in the
product $\prod_{\substack{i\in\left[  n\right]  ;\\i\neq p}}b_{i}$. We thus
obtain%
\begin{equation}
\prod_{\substack{i\in\left[  n\right]  ;\\i\neq p}}b_{i}=\underbrace{\prod
_{\substack{i\in\left[  n\right]  ;\\\tau\left(  i\right)  \neq p}%
}}_{\substack{=\prod_{i=1}^{n-1}\\\text{(by
(\ref{pf.lem.sol.powerdet.gen.taures.short.bprod.preq}))}}}b_{\tau\left(
i\right)  }=\prod_{i=1}^{n-1}b_{\tau\left(  i\right)  }.
\label{pf.lem.sol.powerdet.gen.taures.short.bprod.3}%
\end{equation}
Now,%
\begin{align*}
\prod_{i=1}^{n-1}\underbrace{b_{\left(  \tau\mid_{\left[  n-1\right]
}\right)  \left(  i\right)  }}_{\substack{=b_{\tau\left(  i\right)
}\\\text{(since }\left(  \tau\mid_{\left[  n-1\right]  }\right)  \left(
i\right)  =\tau\left(  i\right)  \text{)}}}  &  =\prod_{i=1}^{n-1}%
b_{\tau\left(  i\right)  }=\prod_{\substack{i\in\left[  n\right]  ;\\i\neq
p}}b_{i}\ \ \ \ \ \ \ \ \ \ \left(  \text{by
(\ref{pf.lem.sol.powerdet.gen.taures.short.bprod.3})}\right) \\
&  =\prod_{\substack{i\in\left\{  1,2,\ldots,n\right\}  ;\\i\neq p}%
}b_{i}\ \ \ \ \ \ \ \ \ \ \left(  \text{since }\left[  n\right]  =\left\{
1,2,\ldots,n\right\}  \right)  .
\end{align*}
This proves (\ref{pf.lem.sol.powerdet.gen.taures.short.bprod}).}.

We have%
\begin{align*}
&  \underbrace{\sum_{\sigma\in S_{n}}}_{\substack{=\sum_{q\in\left[  n\right]
}\sum_{\substack{\sigma\in S_{n};\\\sigma\left(  p\right)  =q}%
}\\\text{(because for each }\sigma\in S_{n}\text{, there}\\\text{exists a
unique }q\in\left[  n\right]  \\\text{such that }\sigma\left(  p\right)
=q\text{)}}}\left(  -1\right)  ^{\sigma}\underbrace{\prod_{i=1}^{n-1}%
a_{\left(  \tau\mid_{\left[  n-1\right]  }\right)  \left(  i\right)
,\sigma\left(  \left(  \tau\mid_{\left[  n-1\right]  }\right)  \left(
i\right)  \right)  }}_{\substack{=\prod_{\substack{i\in\left\{  1,2,\ldots
,n\right\}  ;\\i\neq p}}a_{i,\sigma\left(  i\right)  }\\\text{(by
(\ref{pf.lem.sol.powerdet.gen.taures.short.bprod}), applied to }%
b_{k}=a_{k,\sigma\left(  k\right)  }\text{)}}}\\
&  =\underbrace{\sum_{q\in\left[  n\right]  }}_{=\sum_{q=1}^{n}}%
\underbrace{\sum_{\substack{\sigma\in S_{n};\\\sigma\left(  p\right)
=q}}\left(  -1\right)  ^{\sigma}\prod_{\substack{i\in\left\{  1,2,\ldots
,n\right\}  ;\\i\neq p}}a_{i,\sigma\left(  i\right)  }}_{\substack{=\left(
-1\right)  ^{p+q}\det\left(  A_{\sim p,\sim q}\right)  \\\text{(by Lemma
\ref{lem.laplace.Apq})}}}\\
&  =\sum_{q=1}^{n}\left(  -1\right)  ^{p+q}\det\left(  A_{\sim p,\sim
q}\right)  =\sum_{q=1}^{n}\left(  -1\right)  ^{\tau\left(  n\right)  +q}%
\det\left(  A_{\sim\left(  \tau\left(  n\right)  \right)  ,\sim q}\right)
\end{align*}
(since $p=\tau\left(  n\right)  $). This proves Lemma
\ref{lem.sol.powerdet.gen.taures}.
\end{proof}
\end{vershort}

\begin{verlong}
\begin{proof}
[Proof of Lemma \ref{lem.sol.powerdet.gen.taures}.]We have $\tau\in S_{n}$. In
other words, $\tau$ is a permutation of $\left\{  1,2,\ldots,n\right\}  $
(since $S_{n}$ is the set of all permutations of $\left\{  1,2,\ldots
,n\right\}  $). In other words, $\tau$ is a permutation of $\left[  n\right]
$ (since $\left\{  1,2,\ldots,n\right\}  =\left[  n\right]  $). In other
words, $\tau$ is a bijective map $\left[  n\right]  \rightarrow\left[
n\right]  $.

Also, $n\in\left[  n\right]  $ (since $n\geq1$) and%
\[
\underbrace{\left[  n\right]  }_{=\left\{  1,2,\ldots,n\right\}  }%
\setminus\left\{  n\right\}  =\left\{  1,2,\ldots,n\right\}  \setminus\left\{
n\right\}  =\left\{  1,2,\ldots,n-1\right\}  .
\]
Thus, $\left[  n-1\right]  =\left\{  1,2,\ldots,n-1\right\}  =\left[
n\right]  \setminus\left\{  n\right\}  \subseteq\left[  n\right]  $.

The value $\tau\left(  n\right)  $ is well-defined (since $n\in\left[
n\right]  $). Let $p=\tau\left(  n\right)  $. Then, $p=\tau\left(  n\right)
\in\left[  n\right]  =\left\{  1,2,\ldots,n\right\}  $. Also, $\tau$ is a
bijection (since $\tau$ is a bijective map). Therefore, from $p=\tau\left(
n\right)  $, we obtain $n=\tau^{-1}\left(  p\right)  $.

If $b_{1},b_{2},\ldots,b_{n}$ are $n$ elements of $\mathbb{K}$, then%
\begin{equation}
\prod_{i=1}^{n-1}b_{\left(  \tau\mid_{\left[  n-1\right]  }\right)  \left(
i\right)  }=\prod_{\substack{i\in\left\{  1,2,\ldots,n\right\}  ;\\i\neq
p}}b_{i} \label{pf.lem.sol.powerdet.gen.taures.bprod}%
\end{equation}
\footnote{\textit{Proof of (\ref{pf.lem.sol.powerdet.gen.taures.bprod}):} Let
$b_{1},b_{2},\ldots,b_{n}$ be $n$ elements of $\mathbb{K}$.
\par
For each $i\in\left[  n\right]  $, we have the following chain of logical
equivalences:%
\begin{equation}
\left(  i=\underbrace{n}_{=\tau^{-1}\left(  p\right)  }\right)
\ \Longleftrightarrow\ \left(  i=\tau^{-1}\left(  p\right)  \right)
\ \Longleftrightarrow\ \left(  \tau\left(  i\right)  =p\right)  .
\label{pf.lem.sol.powerdet.gen.taures.bprod.eq1}%
\end{equation}
Thus, for each $i\in\left[  n\right]  $, we have the following chain of
logical equivalences:%
\[
\left(  i\neq n\right)  \ \Longleftrightarrow\ \left(  \text{not
}\underbrace{i=n}_{\substack{\Longleftrightarrow\ \left(  \tau\left(
i\right)  =p\right)  \\\text{(by
(\ref{pf.lem.sol.powerdet.gen.taures.bprod.eq1}))}}}\right)
\ \Longleftrightarrow\ \left(  \text{not }\tau\left(  i\right)  =p\right)
\ \Longleftrightarrow\ \left(  \tau\left(  i\right)  \neq p\right)  .
\]
In other words, for each $i\in\left[  n\right]  $, the condition $\left(
i\neq n\right)  $ holds if and only if the condition $\left(  \tau\left(
i\right)  \neq p\right)  $ holds. Thus, we have the following equality of
product signs:%
\[
\prod_{\substack{i\in\left[  n\right]  ;\\i\neq n}}=\prod_{\substack{i\in
\left[  n\right]  ;\\\tau\left(  i\right)  \neq p}}.
\]
Hence,%
\begin{align}
\prod_{\substack{i\in\left[  n\right]  ;\\\tau\left(  i\right)  \neq p}}  &
=\prod_{\substack{i\in\left[  n\right]  ;\\i\neq n}}=\prod_{i\in\left[
n\right]  \setminus\left\{  n\right\}  }=\prod_{i\in\left\{  1,2,\ldots
,n-1\right\}  }\ \ \ \ \ \ \ \ \ \ \left(  \text{since }\left[  n\right]
\setminus\left\{  n\right\}  =\left\{  1,2,\ldots,n-1\right\}  \right)
\nonumber\\
&  =\prod_{i=1}^{n-1}. \label{pf.lem.sol.powerdet.gen.taures.bprod.preq}%
\end{align}
\par
Now, the map $\tau:\left[  n\right]  \rightarrow\left[  n\right]  $ is a
bijection (since $\tau$ is a bijective map $\left[  n\right]  \rightarrow
\left[  n\right]  $). Thus, we can substitute $\tau\left(  i\right)  $ for $i$
in the product $\prod_{\substack{i\in\left[  n\right]  ;\\i\neq p}}b_{i}$. We
thus obtain%
\begin{equation}
\prod_{\substack{i\in\left[  n\right]  ;\\i\neq p}}b_{i}=\underbrace{\prod
_{\substack{i\in\left[  n\right]  ;\\\tau\left(  i\right)  \neq p}%
}}_{\substack{=\prod_{i=1}^{n-1}\\\text{(by
(\ref{pf.lem.sol.powerdet.gen.taures.bprod.preq}))}}}b_{\tau\left(  i\right)
}=\prod_{i=1}^{n-1}b_{\tau\left(  i\right)  }.
\label{pf.lem.sol.powerdet.gen.taures.bprod.3}%
\end{equation}
Now,%
\begin{align*}
\prod_{i=1}^{n-1}\underbrace{b_{\left(  \tau\mid_{\left[  n-1\right]
}\right)  \left(  i\right)  }}_{\substack{=b_{\tau\left(  i\right)
}\\\text{(since }\left(  \tau\mid_{\left[  n-1\right]  }\right)  \left(
i\right)  =\tau\left(  i\right)  \text{)}}}  &  =\prod_{i=1}^{n-1}%
b_{\tau\left(  i\right)  }=\prod_{\substack{i\in\left[  n\right]  ;\\i\neq
p}}b_{i}\ \ \ \ \ \ \ \ \ \ \left(  \text{by
(\ref{pf.lem.sol.powerdet.gen.taures.bprod.3})}\right) \\
&  =\prod_{\substack{i\in\left\{  1,2,\ldots,n\right\}  ;\\i\neq p}%
}b_{i}\ \ \ \ \ \ \ \ \ \ \left(  \text{since }\left[  n\right]  =\left\{
1,2,\ldots,n\right\}  \right)  .
\end{align*}
This proves (\ref{pf.lem.sol.powerdet.gen.taures.bprod}).}.

We have%
\begin{align*}
&  \underbrace{\sum_{\sigma\in S_{n}}}_{\substack{=\sum_{q\in\left[  n\right]
}\sum_{\substack{\sigma\in S_{n};\\\sigma\left(  p\right)  =q}%
}\\\text{(because for each }\sigma\in S_{n}\text{, there}\\\text{exists a
unique }q\in\left[  n\right]  \\\text{such that }\sigma\left(  p\right)
=q\text{)}}}\left(  -1\right)  ^{\sigma}\underbrace{\prod_{i=1}^{n-1}%
a_{\left(  \tau\mid_{\left[  n-1\right]  }\right)  \left(  i\right)
,\sigma\left(  \left(  \tau\mid_{\left[  n-1\right]  }\right)  \left(
i\right)  \right)  }}_{\substack{=\prod_{\substack{i\in\left\{  1,2,\ldots
,n\right\}  ;\\i\neq p}}a_{i,\sigma\left(  i\right)  }\\\text{(by
(\ref{pf.lem.sol.powerdet.gen.taures.bprod}), applied to }b_{k}=a_{k,\sigma
\left(  k\right)  }\text{)}}}\\
&  =\underbrace{\sum_{q\in\left[  n\right]  }}_{\substack{=\sum_{q\in\left\{
1,2,\ldots,n\right\}  }\\\text{(since }\left[  n\right]  =\left\{
1,2,\ldots,n\right\}  \text{)}}}\sum_{\substack{\sigma\in S_{n};\\\sigma
\left(  p\right)  =q}}\left(  -1\right)  ^{\sigma}\prod_{\substack{i\in
\left\{  1,2,\ldots,n\right\}  ;\\i\neq p}}a_{i,\sigma\left(  i\right)  }\\
&  =\underbrace{\sum_{q\in\left\{  1,2,\ldots,n\right\}  }}_{=\sum_{q=1}^{n}%
}\underbrace{\sum_{\substack{\sigma\in S_{n};\\\sigma\left(  p\right)
=q}}\left(  -1\right)  ^{\sigma}\prod_{\substack{i\in\left\{  1,2,\ldots
,n\right\}  ;\\i\neq p}}a_{i,\sigma\left(  i\right)  }}_{\substack{=\left(
-1\right)  ^{p+q}\det\left(  A_{\sim p,\sim q}\right)  \\\text{(by Lemma
\ref{lem.laplace.Apq})}}}\\
&  =\sum_{q=1}^{n}\left(  -1\right)  ^{p+q}\det\left(  A_{\sim p,\sim
q}\right)  =\sum_{q=1}^{n}\left(  -1\right)  ^{\tau\left(  n\right)  +q}%
\det\left(  A_{\sim\left(  \tau\left(  n\right)  \right)  ,\sim q}\right)
\end{align*}
(since $p=\tau\left(  n\right)  $). This proves Lemma
\ref{lem.sol.powerdet.gen.taures}.
\end{proof}
\end{verlong}

\begin{lemma}
\label{lem.sol.powerdet.gen.n-1}Let $n\geq1$ be an integer.

\textbf{(a)} We have $\left\vert \left\{  \tau\in S_{n}\ \mid\ \tau\left(
n\right)  =n\right\}  \right\vert =\left(  n-1\right)  !$.

\textbf{(b)} Let $p\in\left[  n\right]  $ and $q\in\left[  n\right]  $. Then,%
\[
\left\vert \left\{  \tau\in S_{n}\ \mid\ \tau\left(  p\right)  =q\right\}
\right\vert =\left(  n-1\right)  !.
\]

\end{lemma}

\begin{proof}
[Proof of Lemma \ref{lem.sol.powerdet.gen.n-1}.]We have $n\geq1$. Thus,
$n-1\in\mathbb{N}$.

Hence, Corollary \ref{cor.transpos.code.n!} (applied to $n-1$ instead of $n$)
yields $\left\vert S_{n-1}\right\vert =\left(  n-1\right)  !$.

\textbf{(a)} Clearly, $n$ is a positive integer (since $n$ is an integer
satisfying $n\geq1$).

Define a subset $T$ of $S_{n}$ by $T=\left\{  \tau\in S_{n}\ \mid\ \tau\left(
n\right)  =n\right\}  $. Define a map $\Phi:S_{n-1}\rightarrow T$ as in the
proof of Lemma \ref{lem.laplace.lem}. Then, the map $\Phi$ is a
bijection\footnote{This was shown in the proof of Lemma \ref{lem.laplace.lem}%
.}. Hence, there exists a bijection $S_{n-1}\rightarrow T$ (namely, $\Phi$).
Thus, $\left\vert T\right\vert =\left\vert S_{n-1}\right\vert =\left(
n-1\right)  !$. Since $T=\left\{  \tau\in S_{n}\ \mid\ \tau\left(  n\right)
=n\right\}  $, this rewrites as $\left\vert \left\{  \tau\in S_{n}\ \mid
\ \tau\left(  n\right)  =n\right\}  \right\vert =\left(  n-1\right)  !$. This
proves Lemma \ref{lem.sol.powerdet.gen.n-1} \textbf{(a)}.

\begin{vershort}
\textbf{(b)} Define a map $T:\left\{  \tau\in S_{n}\ \mid\ \tau\left(
n\right)  =n\right\}  \rightarrow\left\{  \tau\in S_{n}\ \mid\ \tau\left(
p\right)  =q\right\}  $ as in Lemma \ref{lem.laplace.gpshort} \textbf{(d)}.
Then, Lemma \ref{lem.laplace.gpshort} \textbf{(d)} shows that this map $T$ is
well-defined and bijective. Thus, there exists a bijection from $\left\{
\tau\in S_{n}\ \mid\ \tau\left(  n\right)  =n\right\}  $ to $\left\{  \tau\in
S_{n}\ \mid\ \tau\left(  p\right)  =q\right\}  $ (namely, $T$). Hence,%
\[
\left\vert \left\{  \tau\in S_{n}\ \mid\ \tau\left(  p\right)  =q\right\}
\right\vert =\left\vert \left\{  \tau\in S_{n}\ \mid\ \tau\left(  n\right)
=n\right\}  \right\vert =\left(  n-1\right)  !
\]
(by Lemma \ref{lem.sol.powerdet.gen.n-1} \textbf{(a)}). This proves Lemma
\ref{lem.sol.powerdet.gen.n-1} \textbf{(b)}. \qedhere

\end{vershort}

\begin{verlong}
\textbf{(b)} Define a map $T:\left\{  \tau\in S_{n}\ \mid\ \tau\left(
n\right)  =n\right\}  \rightarrow\left\{  \tau\in S_{n}\ \mid\ \tau\left(
p\right)  =q\right\}  $ as in Lemma \ref{lem.laplace.gp} \textbf{(g)}. Then,
Lemma \ref{lem.laplace.gp} \textbf{(g)} shows that this map $T$ is
well-defined and bijective. Thus, $T$ is a bijection from $\left\{  \tau\in
S_{n}\ \mid\ \tau\left(  n\right)  =n\right\}  $ to $\left\{  \tau\in
S_{n}\ \mid\ \tau\left(  p\right)  =q\right\}  $. Hence, there exists a
bijection from $\left\{  \tau\in S_{n}\ \mid\ \tau\left(  n\right)
=n\right\}  $ to $\left\{  \tau\in S_{n}\ \mid\ \tau\left(  p\right)
=q\right\}  $ (namely, $T$). Hence,%
\[
\left\vert \left\{  \tau\in S_{n}\ \mid\ \tau\left(  p\right)  =q\right\}
\right\vert =\left\vert \left\{  \tau\in S_{n}\ \mid\ \tau\left(  n\right)
=n\right\}  \right\vert =\left(  n-1\right)  !
\]
(by Lemma \ref{lem.sol.powerdet.gen.n-1} \textbf{(a)}). This proves Lemma
\ref{lem.sol.powerdet.gen.n-1} \textbf{(b)}.
\end{verlong}
\end{proof}

Finally, we are approaching part \textbf{(b)} of Exercise
\ref{exe.powerdet.gen}:

\begin{proposition}
\label{prop.sol.powerdet.gen.b}Let $n\geq1$ be an integer. Let $A=\left(
a_{i,j}\right)  _{1\leq i\leq n,\ 1\leq j\leq n}$ be an $n\times n$-matrix.
Then,%
\[
\sum_{\sigma\in S_{n}}\left(  -1\right)  ^{\sigma}\left(  \sum_{i=1}%
^{n}a_{i,\sigma\left(  i\right)  }\right)  ^{n-1}=\left(  n-1\right)
!\cdot\sum_{p=1}^{n}\sum_{q=1}^{n}\left(  -1\right)  ^{p+q}\det\left(  A_{\sim
p,\sim q}\right)  .
\]

\end{proposition}

\begin{proof}
[Proof of Proposition \ref{prop.sol.powerdet.gen.b}.]From $n\geq1$, we obtain
$n-1\in\mathbb{N}$. Hence, for every $\sigma\in S_{n}$, we have%
\begin{equation}
\left(  \sum_{i=1}^{n}a_{i,\sigma\left(  i\right)  }\right)  ^{n-1}%
=\sum_{\kappa:\left[  n-1\right]  \rightarrow\left[  n\right]  }\prod
_{i=1}^{n-1}a_{\kappa\left(  i\right)  ,\sigma\left(  \kappa\left(  i\right)
\right)  } \label{pf.prop.sol.powerdet.gen.b.prodrule}%
\end{equation}
(by Corollary \ref{cor.sol.powerdet.gen.prodrule}, applied to $k=n-1$).

Lemma \ref{lem.sol.powerdet.gen.kappas} shows that the map%
\begin{align*}
S_{n}  &  \rightarrow\left\{  f:\left[  n-1\right]  \rightarrow\left[
n\right]  \ \mid\ \left\vert f\left(  \left[  n-1\right]  \right)  \right\vert
\geq n-1\right\}  ,\\
\tau &  \mapsto\tau\mid_{\left[  n-1\right]  }%
\end{align*}
is well-defined and bijective. Thus, this map is a bijection.

\begin{verlong}
Notice that $n\in\left[  n\right]  $ (since $n\geq1$).
\end{verlong}

We have $\left[  n\right]  =\left\{  1,2,\ldots,n\right\}  $. Hence, we have
the following equality of summation signs:%
\begin{equation}
\sum_{p\in\left[  n\right]  }=\sum_{p\in\left\{  1,2,\ldots,n\right\}  }%
=\sum_{p=1}^{n}. \label{pf.prop.sol.powerdet.gen.b.sumeqp}%
\end{equation}

Now,%
\begin{align*}
&  \sum_{\sigma\in S_{n}}\left(  -1\right)  ^{\sigma}\underbrace{\left(
\sum_{i=1}^{n}a_{i,\sigma\left(  i\right)  }\right)  ^{n-1}}_{\substack{=\sum
_{\kappa:\left[  n-1\right]  \rightarrow\left[  n\right]  }\prod_{i=1}%
^{n-1}a_{\kappa\left(  i\right)  ,\sigma\left(  \kappa\left(  i\right)
\right)  }\\\text{(by (\ref{pf.prop.sol.powerdet.gen.b.prodrule}))}}}\\
&  =\sum_{\sigma\in S_{n}}\left(  -1\right)  ^{\sigma}\sum_{\kappa:\left[
n-1\right]  \rightarrow\left[  n\right]  }\prod_{i=1}^{n-1}a_{\kappa\left(
i\right)  ,\sigma\left(  \kappa\left(  i\right)  \right)  }=\underbrace{\sum
_{\sigma\in S_{n}}\sum_{\kappa:\left[  n-1\right]  \rightarrow\left[
n\right]  }}_{=\sum_{\kappa:\left[  n-1\right]  \rightarrow\left[  n\right]
}\sum_{\sigma\in S_{n}}}\left(  -1\right)  ^{\sigma}\prod_{i=1}^{n-1}%
a_{\kappa\left(  i\right)  ,\sigma\left(  \kappa\left(  i\right)  \right)  }\\
&  =\sum_{\kappa:\left[  n-1\right]  \rightarrow\left[  n\right]  }%
\sum_{\sigma\in S_{n}}\left(  -1\right)  ^{\sigma}\prod_{i=1}^{n-1}%
a_{\kappa\left(  i\right)  ,\sigma\left(  \kappa\left(  i\right)  \right)  }\\
&  =\sum_{\substack{\kappa:\left[  n-1\right]  \rightarrow\left[  n\right]
;\\\left\vert \kappa\left(  \left[  n-1\right]  \right)  \right\vert
<n-1}}\underbrace{\sum_{\sigma\in S_{n}}\left(  -1\right)  ^{\sigma}%
\prod_{i=1}^{n-1}a_{\kappa\left(  i\right)  ,\sigma\left(  \kappa\left(
i\right)  \right)  }}_{\substack{=0\\\text{(by Lemma
\ref{lem.sol.powerdet.gen.toofew} (applied to }k=n-1\text{))}}}\\
&  \ \ \ \ \ \ \ \ \ \ +\underbrace{\sum_{\substack{\kappa:\left[  n-1\right]
\rightarrow\left[  n\right]  ;\\\left\vert \kappa\left(  \left[  n-1\right]
\right)  \right\vert \geq n-1}}}_{=\sum_{\kappa\in\left\{  f:\left[
n-1\right]  \rightarrow\left[  n\right]  \ \mid\ \left\vert f\left(  \left[
n-1\right]  \right)  \right\vert \geq n-1\right\}  }}\sum_{\sigma\in S_{n}%
}\left(  -1\right)  ^{\sigma}\prod_{i=1}^{n-1}a_{\kappa\left(  i\right)
,\sigma\left(  \kappa\left(  i\right)  \right)  }\\
&  \ \ \ \ \ \ \ \ \ \ \left(
\begin{array}
[c]{c}%
\text{since each map }\kappa:\left[  n-1\right]  \rightarrow\left[  n\right]
\text{ satisfies either }\left\vert \kappa\left(  \left[  n-1\right]  \right)
\right\vert <n-1\\
\text{or }\left\vert \kappa\left(  \left[  n-1\right]  \right)  \right\vert
\geq n-1\text{ (but not both)}%
\end{array}
\right) \\
&  =\underbrace{\sum_{\substack{\kappa:\left[  n-1\right]  \rightarrow\left[
n\right]  ;\\\left\vert \kappa\left(  \left[  n-1\right]  \right)  \right\vert
<n-1}}0}_{=0}+\sum_{\kappa\in\left\{  f:\left[  n-1\right]  \rightarrow\left[
n\right]  \ \mid\ \left\vert f\left(  \left[  n-1\right]  \right)  \right\vert
\geq n-1\right\}  }\sum_{\sigma\in S_{n}}\left(  -1\right)  ^{\sigma}%
\prod_{i=1}^{n-1}a_{\kappa\left(  i\right)  ,\sigma\left(  \kappa\left(
i\right)  \right)  }\\
&  =\sum_{\kappa\in\left\{  f:\left[  n-1\right]  \rightarrow\left[  n\right]
\ \mid\ \left\vert f\left(  \left[  n-1\right]  \right)  \right\vert \geq
n-1\right\}  }\sum_{\sigma\in S_{n}}\left(  -1\right)  ^{\sigma}\prod
_{i=1}^{n-1}a_{\kappa\left(  i\right)  ,\sigma\left(  \kappa\left(  i\right)
\right)  }%
\end{align*}%
\begin{align*}
&  =\sum_{\tau\in S_{n}}\underbrace{\sum_{\sigma\in S_{n}}\left(  -1\right)
^{\sigma}\prod_{i=1}^{n-1}a_{\left(  \tau\mid_{\left[  n-1\right]  }\right)
\left(  i\right)  ,\sigma\left(  \left(  \tau\mid_{\left[  n-1\right]
}\right)  \left(  i\right)  \right)  }}_{\substack{=\sum_{q=1}^{n}\left(
-1\right)  ^{\tau\left(  n\right)  +q}\det\left(  A_{\sim\left(  \tau\left(
n\right)  \right)  ,\sim q}\right)  \\\text{(by Lemma
\ref{lem.sol.powerdet.gen.taures})}}}\\
&  \ \ \ \ \ \ \ \ \ \ \left(
\begin{array}
[c]{c}%
\text{here, we have substituted }\tau\mid_{\left[  n-1\right]  }\text{ for
}\kappa\text{ in the outer sum,}\\
\text{since the}\\
\text{map }S_{n}\rightarrow\left\{  f:\left[  n-1\right]  \rightarrow\left[
n\right]  \ \mid\ \left\vert f\left(  \left[  n-1\right]  \right)  \right\vert
\geq n-1\right\}  ,\ \tau\mapsto\tau\mid_{\left[  n-1\right]  }\\
\text{is a bijection}%
\end{array}
\right) \\
&  =\underbrace{\sum_{\tau\in S_{n}}}_{\substack{=\sum_{p\in\left[  n\right]
}\sum_{\substack{\tau\in S_{n};\\\tau\left(  n\right)  =p}}\\\text{(because
for each }\tau\in S_{n}\text{, there}\\\text{exists a unique }p\in\left[
n\right]  \\\text{such that }\tau\left(  n\right)  =p\text{)}}}\sum_{q=1}%
^{n}\left(  -1\right)  ^{\tau\left(  n\right)  +q}\det\left(  A_{\sim\left(
\tau\left(  n\right)  \right)  ,\sim q}\right) \\
&  =\sum_{p\in\left[  n\right]  }\underbrace{\sum_{\substack{\tau\in
S_{n};\\\tau\left(  n\right)  =p}}\sum_{q=1}^{n}}_{=\sum_{q=1}^{n}%
\sum_{\substack{\tau\in S_{n};\\\tau\left(  n\right)  =p}}}\underbrace{\left(
-1\right)  ^{\tau\left(  n\right)  +q}\det\left(  A_{\sim\left(  \tau\left(
n\right)  \right)  ,\sim q}\right)  }_{\substack{=\left(  -1\right)
^{p+q}\det\left(  A_{\sim p,\sim q}\right)  \\\text{(since }\tau\left(
n\right)  =p\text{)}}}\\
&  =\sum_{p\in\left[  n\right]  }\sum_{q=1}^{n}\underbrace{\sum
_{\substack{\tau\in S_{n};\\\tau\left(  n\right)  =p}}\left(  -1\right)
^{p+q}\det\left(  A_{\sim p,\sim q}\right)  }_{=\left\vert \left\{  \tau\in
S_{n}\ \mid\ \tau\left(  n\right)  =p\right\}  \right\vert \left(  -1\right)
^{p+q}\det\left(  A_{\sim p,\sim q}\right)  }\\
&  =\underbrace{\sum_{p\in\left[  n\right]  }}_{\substack{=\sum_{p=1}%
^{n}\\\text{(by (\ref{pf.prop.sol.powerdet.gen.b.sumeqp}))}}}\sum_{q=1}%
^{n}\underbrace{\left\vert \left\{  \tau\in S_{n}\ \mid\ \tau\left(  n\right)
=p\right\}  \right\vert }_{\substack{=\left(  n-1\right)  !\\\text{(by Lemma
\ref{lem.sol.powerdet.gen.n-1} \textbf{(b)}}\\\text{(applied to }n\text{ and
}p\text{ instead of }p\text{ and }q\text{))}}}\left(  -1\right)  ^{p+q}%
\det\left(  A_{\sim p,\sim q}\right) \\
&  =\sum_{p=1}^{n}\sum_{q=1}^{n}\left(  n-1\right)  !\left(  -1\right)
^{p+q}\det\left(  A_{\sim p,\sim q}\right)  =\left(  n-1\right)  !\cdot
\sum_{p=1}^{n}\sum_{q=1}^{n}\left(  -1\right)  ^{p+q}\det\left(  A_{\sim
p,\sim q}\right)  .
\end{align*}
This proves Proposition \ref{prop.sol.powerdet.gen.b}.
\end{proof}

\begin{proof}
[Solution to Exercise \ref{exe.powerdet.gen}.]Exercise \ref{exe.powerdet.gen}
\textbf{(a)} follows from Proposition \ref{prop.sol.powerdet.gen.a}.

Exercise \ref{exe.powerdet.gen} \textbf{(b)} follows from Proposition
\ref{prop.sol.powerdet.gen.b}.
\end{proof}

\subsection{Solution to Exercise \ref{exe.powerdet.r1}}

\subsubsection{Solving the exercise}

We prepare for the solution to Exercise \ref{exe.powerdet.r1} by showing a few lemmas.

\begin{lemma}
\label{lem.sol.powerdet.r1.EI}Let $n\in\mathbb{N}$. We let $\mathbf{E}$ be the
subset%
\[
\left\{  \left(  k_{1},k_{2},\ldots,k_{n}\right)  \in\mathbb{N}^{n}%
\ \mid\ \text{the integers }k_{1},k_{2},\ldots,k_{n}\text{ are distinct}%
\right\}
\]
of $\mathbb{N}^{n}$. We let $\mathbf{I}$ be the subset%
\[
\left\{  \left(  k_{1},k_{2},\ldots,k_{n}\right)  \in\mathbb{N}^{n}%
\ \mid\ k_{1}<k_{2}<\cdots<k_{n}\right\}
\]
of $\mathbb{N}^{n}$. Then, the map%
\begin{align*}
\mathbf{I}\times S_{n}  &  \rightarrow\mathbf{E},\\
\left(  \left(  g_{1},g_{2},\ldots,g_{n}\right)  ,\sigma\right)   &
\mapsto\left(  g_{\sigma\left(  1\right)  },g_{\sigma\left(  2\right)
},\ldots,g_{\sigma\left(  n\right)  }\right)
\end{align*}
is well-defined and is a bijection.
\end{lemma}

Lemma \ref{lem.sol.powerdet.r1.EI} is identical to Lemma
\ref{lem.cauchy-binet.EI} except that it uses the set $\mathbb{N}$ instead of
the set $\left[  m\right]  $ used in the latter lemma. Thus, it should not be
a surprise that its proof is the same as the proof of the latter lemma:

\begin{proof}
[Proof of Lemma \ref{lem.sol.powerdet.r1.EI}.]In order to obtain a proof of
Lemma \ref{lem.sol.powerdet.r1.EI}, it is sufficient to replace every
appearance of \textquotedblleft$\left[  m\right]  $\textquotedblright\ by
\textquotedblleft$\mathbb{N}$\textquotedblright\ in the proof of Lemma
\ref{lem.cauchy-binet.EI}. (Of course, the notation $\left[  n\right]  $
should be understood to mean the set $\left\{  1,2,\ldots,n\right\}  $ in this proof.)
\end{proof}

\begin{lemma}
\label{lem.sol.powerdet.r1.0}Let $n\in\mathbb{N}$. Let $\left(  k_{1}%
,k_{2},\ldots,k_{n}\right)  \in\mathbb{N}^{n}$ be such that $k_{1}%
<k_{2}<\cdots<k_{n}$.

\textbf{(a)} We have $k_{1}+k_{2}+\cdots+k_{n}\geq\dbinom{n}{2}$.

\textbf{(b)} If $k_{1}+k_{2}+\cdots+k_{n}\leq\dbinom{n}{2}$, then $\left(
k_{1},k_{2},\ldots,k_{n}\right)  =\left(  0,1,\ldots,n-1\right)  $.
\end{lemma}

\begin{vershort}
\begin{proof}
[Proof of Lemma \ref{lem.sol.powerdet.r1.0}.]\textbf{(b)} Assume that
$k_{1}+k_{2}+\cdots+k_{n}\leq\dbinom{n}{2}$. We must prove that $\left(
k_{1},k_{2},\ldots,k_{n}\right)  =\left(  0,1,\ldots,n-1\right)  $.

For each $i\in\left\{  1,2,\ldots,n\right\}  $, define an integer $z_{i}$ by
$z_{i}=k_{i}-i$.

For each $i\in\left\{  1,2,\ldots,n-1\right\}  $, we have $z_{i}\leq z_{i+1}%
$\ \ \ \ \footnote{\textit{Proof.} Let $i\in\left\{  1,2,\ldots,n-1\right\}
$. The definition of $z_{i}$ yields $z_{i}=k_{i}-i$. The definition of
$z_{i+1}$ yields $z_{i+1}=k_{i+1}-\left(  i+1\right)  $. But $k_{i}<k_{i+1}$
(since $k_{1}<k_{2}<\cdots<k_{n}$), and thus $k_{i}\leq k_{i+1}-1$ (since
$k_{i}$ and $k_{i+1}$ are integers). Now, $z_{i}=\underbrace{k_{i}}_{\leq
k_{i+1}-1}-i\leq k_{i+1}-1-i=k_{i+1}-\left(  i+1\right)  =z_{i+1}$. Qed.}. In
other words, $z_{1}\leq z_{2}\leq\cdots\leq z_{n}$.

For each $j\in\left\{  1,2,\ldots,n\right\}  $, we have%
\begin{equation}
k_{j}\geq j-1 \label{pf.lem.sol.powerdet.r1.0.short.2}%
\end{equation}
\footnote{\textit{Proof of (\ref{pf.lem.sol.powerdet.r1.0.short.2}):} Let
$j\in\left\{  1,2,\ldots,n\right\}  $. Thus, $1\leq j\leq n$. Hence, $n\geq1$,
so that $1\in\left\{  1,2,\ldots,n\right\}  $.
\par
But recall that $z_{1}\leq z_{2}\leq\cdots\leq z_{n}$. In other words, we have
$z_{u}\leq z_{v}$ whenever $u$ and $v$ are two elements of $\left\{
1,2,\ldots,n\right\}  $ satisfying $u\leq v$. Applying this to $u=1$ and
$v=j$, we obtain $z_{1}\leq z_{j}$ (since $1\leq j$). But the definition of
$z_{1}$ yields $z_{1}=\underbrace{k_{1}}_{\substack{\geq0\\\text{(since }%
k_{1}\in\mathbb{N}\text{)}}}-1\geq0-1=-1$. Hence, $-1\leq z_{1}\leq
z_{j}=k_{j}-j$ (by the definition of $z_{j}$). Hence, $k_{j}-j\geq-1$, so that
$k_{j}\geq j-1$. This proves (\ref{pf.lem.sol.powerdet.r1.0.short.2}).}.

Now, each $j\in\left\{  1,2,\ldots,n\right\}  $ satisfies%
\begin{equation}
k_{j}=j-1. \label{pf.lem.sol.powerdet.r1.0.short.3}%
\end{equation}

[\textit{Proof of (\ref{pf.lem.sol.powerdet.r1.0.short.3}):} Let $j\in\left\{
1,2,\ldots,n\right\}  $. We must prove that $k_{j}=j-1$.

Assume the contrary. Thus, $k_{j}\neq j-1$. Combining this with
(\ref{pf.lem.sol.powerdet.r1.0.short.2}), we obtain $k_{j}>j-1$.

Now, $k_{1}+k_{2}+\cdots+k_{n}\leq\dbinom{n}{2}$, so that%
\begin{align*}
\dbinom{n}{2}  &  \geq k_{1}+k_{2}+\cdots+k_{n}=\sum_{i\in\left\{
1,2,\ldots,n\right\}  }k_{i}\\
&  =\underbrace{k_{j}}_{>j-1}+\sum_{\substack{i\in\left\{  1,2,\ldots
,n\right\}  ;\\i\neq j}}\underbrace{k_{i}}_{\substack{\geq i-1\\\text{(by
(\ref{pf.lem.sol.powerdet.r1.0.short.2}) (applied}\\\text{to }i\text{ instead
of }j\text{))}}}\ \ \ \ \ \ \ \ \ \ \left(
\begin{array}
[c]{c}%
\text{here, we have split off the addend}\\
\text{for }i=j\text{ from the sum}%
\end{array}
\right) \\
&  >\left(  j-1\right)  +\sum_{\substack{i\in\left\{  1,2,\ldots,n\right\}
;\\i\neq j}}\left(  i-1\right)  =\sum_{i\in\left\{  1,2,\ldots,n\right\}
}\left(  i-1\right) \\
&  =0+1+\cdots+\left(  n-1\right)  =\sum_{r=0}^{n-1}r=\dbinom{n}{2}%
\end{align*}
(by Lemma \ref{lem.sol.vander-hook.gauss}). This is absurd. Thus, we have
found a contradiction. This contradiction shows that our assumption was wrong.
Hence, $k_{j}=j-1$ is proven. This proves
(\ref{pf.lem.sol.powerdet.r1.0.short.3}).]

Now, (\ref{pf.lem.sol.powerdet.r1.0.short.3}) shows that $\left(  k_{1}%
,k_{2},\ldots,k_{n}\right)  =\left(  1-1,2-1,\ldots,n-1\right)  =\left(
0,1,\ldots,n-1\right)  $. This proves Lemma \ref{lem.sol.powerdet.r1.0}
\textbf{(b)}.

\textbf{(a)} Assume the contrary. Thus, $k_{1}+k_{2}+\cdots+k_{n}<\dbinom
{n}{2}$. Hence, $k_{1}+k_{2}+\cdots+k_{n}\leq\dbinom{n}{2}$. Therefore, Lemma
\ref{lem.sol.powerdet.r1.0} \textbf{(b)} shows that $\left(  k_{1}%
,k_{2},\ldots,k_{n}\right)  =\left(  0,1,\ldots,n-1\right)  $. Thus,%
\[
k_{1}+k_{2}+\cdots+k_{n}=0+1+\cdots+\left(  n-1\right)  =\sum_{r=0}%
^{n-1}r=\dbinom{n}{2}%
\]
(by Lemma \ref{lem.sol.vander-hook.gauss}), which contradicts $k_{1}%
+k_{2}+\cdots+k_{n}<\dbinom{n}{2}$. This contradiction shows that our
assumption was wrong. Hence, Lemma \ref{lem.sol.powerdet.r1.0} \textbf{(a)} is proven.
\end{proof}
\end{vershort}

\begin{verlong}
\begin{proof}
[Proof of Lemma \ref{lem.sol.powerdet.r1.0}.]We have $\left(  k_{1}%
,k_{2},\ldots,k_{n}\right)  \in\mathbb{N}^{n}$. Thus, $k_{1},k_{2}%
,\ldots,k_{n}$ are $n$ elements of $\mathbb{N}$, therefore $n$ integers.

We observe that
\begin{equation}
0+1+\cdots+\left(  n-1\right)  =\sum_{r=0}^{n-1}r=\dbinom{n}{2}
\label{pf.lem.sol.powerdet.r1.0.sum}%
\end{equation}
(by Lemma \ref{lem.sol.vander-hook.gauss}).

\textbf{(b)} Assume that $k_{1}+k_{2}+\cdots+k_{n}\leq\dbinom{n}{2}$. We must
prove that $\left(  k_{1},k_{2},\ldots,k_{n}\right)  =\left(  0,1,\ldots
,n-1\right)  $.

For each $i\in\left\{  1,2,\ldots,n\right\}  $, define an integer $z_{i}$ by
$z_{i}=k_{i}-i$.

For each $i\in\left\{  1,2,\ldots,n-1\right\}  $, we have $z_{i}\leq z_{i+1}%
$\ \ \ \ \footnote{\textit{Proof.} Let $i\in\left\{  1,2,\ldots,n-1\right\}
$. Then, both $i$ and $i+1$ are elements of $\left\{  1,2,\ldots,n\right\}  $.
Hence, the integers $z_{i}$ and $z_{i+1}$ are well-defined.
\par
The definition of $z_{i}$ yields $z_{i}=k_{i}-i$. The definition of $z_{i+1}$
yields $z_{i+1}=k_{i+1}-\left(  i+1\right)  $. But $k_{i}<k_{i+1}$ (since
$k_{1}<k_{2}<\cdots<k_{n}$), and thus $k_{i}\leq k_{i+1}-1$ (since $k_{i}$ and
$k_{i+1}$ are integers). Now, $z_{i}=\underbrace{k_{i}}_{\leq k_{i+1}-1}-i\leq
k_{i+1}-1-i=k_{i+1}-\left(  i+1\right)  =z_{i+1}$. Qed.}. In other words,
$z_{1}\leq z_{2}\leq\cdots\leq z_{n}$. In other words,%
\begin{equation}
z_{i}\leq z_{j} \label{pf.lem.sol.powerdet.r1.0.1}%
\end{equation}
for every $i\in\left\{  1,2,\ldots,n\right\}  $ and $j\in\left\{
1,2,\ldots,n\right\}  $ satisfying $i\leq j$.

For each $j\in\left\{  1,2,\ldots,n\right\}  $, we have%
\begin{equation}
k_{j}\geq j-1 \label{pf.lem.sol.powerdet.r1.0.2}%
\end{equation}
\footnote{\textit{Proof of (\ref{pf.lem.sol.powerdet.r1.0.2}):} Let
$j\in\left\{  1,2,\ldots,n\right\}  $. Thus, $1\leq j\leq n$. Hence, $n\geq1$,
so that $1\in\left\{  1,2,\ldots,n\right\}  $. Moreover, $1\leq j$. Thus,
(\ref{pf.lem.sol.powerdet.r1.0.1}) (applied to $i=1$) yields $z_{1}\leq z_{j}%
$. But the definition of $z_{1}$ yields $z_{1}=\underbrace{k_{1}%
}_{\substack{\geq0\\\text{(since }k_{1}\in\mathbb{N}\text{)}}}-1\geq0-1=-1$.
Hence, $-1\leq z_{1}\leq z_{j}=k_{j}-j$ (by the definition of $z_{j}$). Hence,
$k_{j}-j\geq-1$, so that $k_{j}\geq j-1$. This proves
(\ref{pf.lem.sol.powerdet.r1.0.2}).}.

Now, each $j\in\left\{  1,2,\ldots,n\right\}  $ satisfies%
\begin{equation}
k_{j}=j-1. \label{pf.lem.sol.powerdet.r1.0.3}%
\end{equation}

[\textit{Proof of (\ref{pf.lem.sol.powerdet.r1.0.3}):} Let $j\in\left\{
1,2,\ldots,n\right\}  $. We must prove that $k_{j}=j-1$.

Assume the contrary. Thus, $k_{j}\neq j-1$. But $k_{j}\geq j-1$ (by
(\ref{pf.lem.sol.powerdet.r1.0.2})). Combining this with $k_{j}\neq j-1$, we
obtain $k_{j}>j-1$.

Now, $k_{1}+k_{2}+\cdots+k_{n}\leq\dbinom{n}{2}$, so that%
\begin{align}
\dbinom{n}{2}  &  \geq k_{1}+k_{2}+\cdots+k_{n}=\underbrace{\sum_{i=1}^{n}%
}_{=\sum_{i\in\left\{  1,2,\ldots,n\right\}  }}k_{i}=\sum_{i\in\left\{
1,2,\ldots,n\right\}  }k_{i}\nonumber\\
&  =\underbrace{k_{j}}_{>j-1}+\sum_{\substack{i\in\left\{  1,2,\ldots
,n\right\}  ;\\i\neq j}}\underbrace{k_{i}}_{\substack{\geq i-1\\\text{(by
(\ref{pf.lem.sol.powerdet.r1.0.2}) (applied}\\\text{to }i\text{ instead of
}j\text{))}}}\ \ \ \ \ \ \ \ \ \ \left(
\begin{array}
[c]{c}%
\text{here, we have split off the addend}\\
\text{for }i=j\text{ from the sum}%
\end{array}
\right) \nonumber\\
&  >\left(  j-1\right)  +\sum_{\substack{i\in\left\{  1,2,\ldots,n\right\}
;\\i\neq j}}\left(  i-1\right)  . \label{pf.lem.sol.powerdet.r1.0.3.pf.1}%
\end{align}
But
\begin{align*}
\dbinom{n}{2}  &  =0+1+\cdots+\left(  n-1\right)  \ \ \ \ \ \ \ \ \ \ \left(
\text{by (\ref{pf.lem.sol.powerdet.r1.0.sum})}\right) \\
&  =\left(  1-1\right)  +\left(  2-1\right)  +\cdots+\left(  n-1\right)
=\underbrace{\sum_{i=1}^{n}}_{=\sum_{i\in\left\{  1,2,\ldots,n\right\}  }%
}\left(  i-1\right)  =\sum_{i\in\left\{  1,2,\ldots,n\right\}  }\left(
i-1\right) \\
&  =\left(  j-1\right)  +\sum_{\substack{i\in\left\{  1,2,\ldots,n\right\}
;\\i\neq j}}\left(  i-1\right)  \ \ \ \ \ \ \ \ \ \ \left(
\begin{array}
[c]{c}%
\text{here, we have split off the addend}\\
\text{for }i=j\text{ from the sum}%
\end{array}
\right) \\
&  <\dbinom{n}{2}\ \ \ \ \ \ \ \ \ \ \left(  \text{by
(\ref{pf.lem.sol.powerdet.r1.0.3.pf.1})}\right)  .
\end{align*}
This is absurd. Thus, we have found a contradiction. This contradiction shows
that our assumption was wrong. Hence, $k_{j}=j-1$ is proven. This proves
(\ref{pf.lem.sol.powerdet.r1.0.3}).]

Now, (\ref{pf.lem.sol.powerdet.r1.0.3}) shows that $\left(  k_{1},k_{2}%
,\ldots,k_{n}\right)  =\left(  1-1,2-1,\ldots,n-1\right)  =\left(
0,1,\ldots,n-1\right)  $. This proves Lemma \ref{lem.sol.powerdet.r1.0}
\textbf{(b)}.

\textbf{(a)} Assume the contrary. Thus, $k_{1}+k_{2}+\cdots+k_{n}<\dbinom
{n}{2}$. Hence, $k_{1}+k_{2}+\cdots+k_{n}\leq\dbinom{n}{2}$. Therefore, Lemma
\ref{lem.sol.powerdet.r1.0} \textbf{(b)} shows that $\left(  k_{1}%
,k_{2},\ldots,k_{n}\right)  =\left(  0,1,\ldots,n-1\right)  $. Thus,%
\[
k_{1}+k_{2}+\cdots+k_{n}=0+1+\cdots+\left(  n-1\right)  =\dbinom{n}%
{2}\ \ \ \ \ \ \ \ \ \ \left(  \text{by (\ref{pf.lem.sol.powerdet.r1.0.sum}%
)}\right)  ,
\]
so that $\dbinom{n}{2}=k_{1}+k_{2}+\cdots+k_{n}<\dbinom{n}{2}$. This is
absurd. Thus, we have found a contradiction. This contradiction shows that our
assumption was wrong. Hence, Lemma \ref{lem.sol.powerdet.r1.0} \textbf{(a)} is proven.
\end{proof}
\end{verlong}

\begin{lemma}
\label{lem.sol.powerdet.r1.1}Let $n\in\mathbb{N}$. Let $a_{1},a_{2}%
,\ldots,a_{n}$ be $n$ elements of $\mathbb{K}$. Let $b_{1},b_{2},\ldots,b_{n}$
be $n$ elements of $\mathbb{K}$.

Let $\left(  k_{1},k_{2},\ldots,k_{n}\right)  \in\mathbb{N}^{n}$.

\textbf{(a)} We have
\[
\sum_{\sigma\in S_{n}}\left(  -1\right)  ^{\sigma}\prod_{i=1}^{n}\left(
a_{i}b_{\sigma\left(  i\right)  }\right)  ^{k_{i}}=\left(  \prod_{i=1}%
^{n}a_{i}^{k_{i}}\right)  \det\left(  \left(  b_{i}^{k_{j}}\right)  _{1\leq
i\leq n,\ 1\leq j\leq n}\right)  .
\]

\textbf{(b)} If the integers $k_{1},k_{2},\ldots,k_{n}$ are not distinct, then%
\[
\sum_{\sigma\in S_{n}}\left(  -1\right)  ^{\sigma}\prod_{i=1}^{n}\left(
a_{i}b_{\sigma\left(  i\right)  }\right)  ^{k_{i}}=0.
\]

\end{lemma}

\begin{proof}
[Proof of Lemma \ref{lem.sol.powerdet.r1.1}.]\textbf{(a)} Define an $n\times
n$-matrix $A$ by $A=\left(  b_{j}^{k_{i}}\right)  _{1\leq i\leq n,\ 1\leq
j\leq n}$. Then, the equality (\ref{eq.det.eq.2}) (applied to $b_{j}^{k_{i}}$
instead of $a_{i,j}$) yields%
\[
\det A=\sum_{\sigma\in S_{n}}\left(  -1\right)  ^{\sigma}\prod_{i=1}%
^{n}b_{\sigma\left(  i\right)  }^{k_{i}}.
\]
But $A=\left(  b_{j}^{k_{i}}\right)  _{1\leq i\leq n,\ 1\leq j\leq n}$. Thus,
$A^{T}=\left(  b_{i}^{k_{j}}\right)  _{1\leq i\leq n,\ 1\leq j\leq n}$ (by the
definition of the transpose matrix $A^{T}$). Exercise \ref{exe.ps4.4} now
yields
\begin{equation}
\det\left(  A^{T}\right)  =\det A=\sum_{\sigma\in S_{n}}\left(  -1\right)
^{\sigma}\prod_{i=1}^{n}b_{\sigma\left(  i\right)  }^{k_{i}}.
\label{pf.lem.sol.powerdet.r1.1.1}%
\end{equation}

Now,%
\begin{align*}
&  \sum_{\sigma\in S_{n}}\left(  -1\right)  ^{\sigma}\prod_{i=1}%
^{n}\underbrace{\left(  a_{i}b_{\sigma\left(  i\right)  }\right)  ^{k_{i}}%
}_{=a_{i}^{k_{i}}b_{\sigma\left(  i\right)  }^{k_{i}}}\\
&  =\sum_{\sigma\in S_{n}}\left(  -1\right)  ^{\sigma}\underbrace{\prod
_{i=1}^{n}\left(  a_{i}^{k_{i}}b_{\sigma\left(  i\right)  }^{k_{i}}\right)
}_{=\left(  \prod_{i=1}^{n}a_{i}^{k_{i}}\right)  \left(  \prod_{i=1}%
^{n}b_{\sigma\left(  i\right)  }^{k_{i}}\right)  }=\sum_{\sigma\in S_{n}%
}\left(  -1\right)  ^{\sigma}\left(  \prod_{i=1}^{n}a_{i}^{k_{i}}\right)
\left(  \prod_{i=1}^{n}b_{\sigma\left(  i\right)  }^{k_{i}}\right) \\
&  =\left(  \prod_{i=1}^{n}a_{i}^{k_{i}}\right)  \underbrace{\sum_{\sigma\in
S_{n}}\left(  -1\right)  ^{\sigma}\prod_{i=1}^{n}b_{\sigma\left(  i\right)
}^{k_{i}}}_{\substack{=\det\left(  A^{T}\right)  \\\text{(by
(\ref{pf.lem.sol.powerdet.r1.1.1}))}}}=\left(  \prod_{i=1}^{n}a_{i}^{k_{i}%
}\right)  \det\left(  \underbrace{A^{T}}_{=\left(  b_{i}^{k_{j}}\right)
_{1\leq i\leq n,\ 1\leq j\leq n}}\right) \\
&  =\left(  \prod_{i=1}^{n}a_{i}^{k_{i}}\right)  \det\left(  \left(
b_{i}^{k_{j}}\right)  _{1\leq i\leq n,\ 1\leq j\leq n}\right)  .
\end{align*}
This proves Lemma \ref{lem.sol.powerdet.r1.1} \textbf{(a)}.

\begin{vershort}
\textbf{(b)} Assume that the integers $k_{1},k_{2},\ldots,k_{n}$ are not
distinct. Thus, there exist two distinct elements $p$ and $q$ of $\left\{
1,2,\ldots,n\right\}  $ such that $k_{p}=k_{q}$. Consider these $p$ and $q$.
Now,%
\begin{align*}
&  \left(  \text{the }p\text{-th column of the matrix }\left(  b_{i}^{k_{j}%
}\right)  _{1\leq i\leq n,\ 1\leq j\leq n}\right) \\
&  =\left(  b_{i}^{k_{p}}\right)  _{1\leq i\leq n,\ 1\leq j\leq1}=\left(
b_{i}^{k_{q}}\right)  _{1\leq i\leq n,\ 1\leq j\leq n}%
\ \ \ \ \ \ \ \ \ \ \left(  \text{since }k_{p}=k_{q}\right) \\
&  =\left(  \text{the }q\text{-th column of the matrix }\left(  b_{i}^{k_{j}%
}\right)  _{1\leq i\leq n,\ 1\leq j\leq n}\right)  .
\end{align*}
Thus, the matrix $\left(  b_{i}^{k_{j}}\right)  _{1\leq i\leq n,\ 1\leq j\leq
n}$ has two equal columns (namely, the $p$-th column and the $q$-th column),
because $p$ and $q$ are distinct. Therefore, Exercise \ref{exe.ps4.6}
\textbf{(f)} (applied to $\left(  b_{i}^{k_{j}}\right)  _{1\leq i\leq
n,\ 1\leq j\leq n}$ instead of $A$) yields $\det\left(  \left(  b_{i}^{k_{j}%
}\right)  _{1\leq i\leq n,\ 1\leq j\leq n}\right)  =0$.

Now, Lemma \ref{lem.sol.powerdet.r1.1} \textbf{(a)} yields%
\[
\sum_{\sigma\in S_{n}}\left(  -1\right)  ^{\sigma}\prod_{i=1}^{n}\left(
a_{i}b_{\sigma\left(  i\right)  }\right)  ^{k_{i}}=\left(  \prod_{i=1}%
^{n}a_{i}^{k_{i}}\right)  \underbrace{\det\left(  \left(  b_{i}^{k_{j}%
}\right)  _{1\leq i\leq n,\ 1\leq j\leq n}\right)  }_{=0}=0.
\]
This proves Lemma \ref{lem.sol.powerdet.r1.1} \textbf{(b)}. \qedhere

\end{vershort}

\begin{verlong}
\textbf{(b)} Assume that the integers $k_{1},k_{2},\ldots,k_{n}$ are not
distinct. Thus, there exist two distinct elements $p$ and $q$ of $\left\{
1,2,\ldots,n\right\}  $ such that $k_{p}=k_{q}$%
\ \ \ \ \footnote{\textit{Proof.} Assume the contrary. Thus, there do not
exist two distinct elements $p$ and $q$ of $\left\{  1,2,\ldots,n\right\}  $
such that $k_{p}=k_{q}$. In other words, every two distinct elements $p$ and
$q$ of $\left\{  1,2,\ldots,n\right\}  $ satisfy $k_{p}\neq k_{q}$. In other
words, the integers $k_{1},k_{2},\ldots,k_{n}$ are distinct. This contradicts
the fact that the integers $k_{1},k_{2},\ldots,k_{n}$ are not distinct. This
contradiction proves that our assumption was wrong. Qed.}. Consider these $p$
and $q$. Now,%
\begin{align*}
&  \left(  \text{the }p\text{-th column of the matrix }\left(  b_{i}^{k_{j}%
}\right)  _{1\leq i\leq n,\ 1\leq j\leq n}\right) \\
&  =\left(  b_{i}^{k_{p}}\right)  _{1\leq i\leq n,\ 1\leq j\leq1}=\left(
b_{i}^{k_{q}}\right)  _{1\leq i\leq n,\ 1\leq j\leq n}%
\end{align*}
(since $k_{p}=k_{q}$). Comparing this with%
\[
\left(  \text{the }q\text{-th column of the matrix }\left(  b_{i}^{k_{j}%
}\right)  _{1\leq i\leq n,\ 1\leq j\leq n}\right)  =\left(  b_{i}^{k_{q}%
}\right)  _{1\leq i\leq n,\ 1\leq j\leq n},
\]
we conclude that%
\begin{align*}
&  \left(  \text{the }p\text{-th column of the matrix }\left(  b_{i}^{k_{j}%
}\right)  _{1\leq i\leq n,\ 1\leq j\leq n}\right) \\
&  =\left(  \text{the }q\text{-th column of the matrix }\left(  b_{i}^{k_{j}%
}\right)  _{1\leq i\leq n,\ 1\leq j\leq n}\right)  .
\end{align*}
Thus, the matrix $\left(  b_{i}^{k_{j}}\right)  _{1\leq i\leq n,\ 1\leq j\leq
n}$ has two equal columns (namely, the $p$-th column and the $q$-th column),
because $p$ and $q$ are distinct. Therefore, Exercise \ref{exe.ps4.6}
\textbf{(f)} (applied to $\left(  b_{i}^{k_{j}}\right)  _{1\leq i\leq
n,\ 1\leq j\leq n}$ instead of $A$) yields $\det\left(  \left(  b_{i}^{k_{j}%
}\right)  _{1\leq i\leq n,\ 1\leq j\leq n}\right)  =0$.

Now, Lemma \ref{lem.sol.powerdet.r1.1} \textbf{(a)} yields%
\[
\sum_{\sigma\in S_{n}}\left(  -1\right)  ^{\sigma}\prod_{i=1}^{n}\left(
a_{i}b_{\sigma\left(  i\right)  }\right)  ^{k_{i}}=\left(  \prod_{i=1}%
^{n}a_{i}^{k_{i}}\right)  \underbrace{\det\left(  \left(  b_{i}^{k_{j}%
}\right)  _{1\leq i\leq n,\ 1\leq j\leq n}\right)  }_{=0}=0.
\]
This proves Lemma \ref{lem.sol.powerdet.r1.1} \textbf{(b)}.
\end{verlong}
\end{proof}

We shall use the notation introduced in Exercise \ref{exe.multinom2} from now
on until the end of the present section.

\begin{lemma}
\label{lem.sol.powerdet.r1.perm-m}Let $n\in\mathbb{N}$. Let $\left(
g_{1},g_{2},\ldots,g_{n}\right)  \in\mathbb{N}^{n}$ and $\sigma\in S_{n}$.

\textbf{(a)} We have $g_{\sigma\left(  1\right)  }+g_{\sigma\left(  2\right)
}+\cdots+g_{\sigma\left(  n\right)  }=g_{1}+g_{2}+\cdots+g_{n}$.

\textbf{(b)} We have $\mathbf{m}\left(  g_{\sigma\left(  1\right)  }%
,g_{\sigma\left(  2\right)  },\ldots,g_{\sigma\left(  n\right)  }\right)
=\mathbf{m}\left(  g_{1},g_{2},\ldots,g_{n}\right)  $.
\end{lemma}

\begin{vershort}
\begin{proof}
[Proof of Lemma \ref{lem.sol.powerdet.r1.perm-m}.]We have $\sigma\in S_{n}$.
In other words, $\sigma$ is a permutation of the set $\left\{  1,2,\ldots
,n\right\}  $ (since $S_{n}$ is the set of all permutations of the set
$\left\{  1,2,\ldots,n\right\}  $). In other words, $\sigma$ is a bijection
$\left\{  1,2,\ldots,n\right\}  \rightarrow\left\{  1,2,\ldots,n\right\}  $.
Hence, we can substitute $\sigma\left(  i\right)  $ for $i$ in the sum
$\sum_{i\in\left\{  1,2,\ldots,n\right\}  }g_{i}$. We thus obtain%
\[
\sum_{i\in\left\{  1,2,\ldots,n\right\}  }g_{i}=\sum_{i\in\left\{
1,2,\ldots,n\right\}  }g_{\sigma\left(  i\right)  }=g_{\sigma\left(  1\right)
}+g_{\sigma\left(  2\right)  }+\cdots+g_{\sigma\left(  n\right)  }.
\]
Hence,%
\begin{equation}
g_{\sigma\left(  1\right)  }+g_{\sigma\left(  2\right)  }+\cdots
+g_{\sigma\left(  n\right)  }=\sum_{i\in\left\{  1,2,\ldots,n\right\}  }%
g_{i}=g_{1}+g_{2}+\cdots+g_{n}. \label{pf.lem.sol.powerdet.r1.perm-m.short.1}%
\end{equation}
This proves Lemma \ref{lem.sol.powerdet.r1.perm-m} \textbf{(a)}.

\textbf{(b)} We have proven that $g_{\sigma\left(  1\right)  }+g_{\sigma
\left(  2\right)  }+\cdots+g_{\sigma\left(  n\right)  }=g_{1}+g_{2}%
+\cdots+g_{n}$. A similar argument shows that
\begin{equation}
g_{\sigma\left(  1\right)  }!g_{\sigma\left(  2\right)  }!\cdots
g_{\sigma\left(  n\right)  }!=g_{1}!g_{2}!\cdots g_{n}!.
\label{pf.lem.sol.powerdet.r1.perm-m.short.2}%
\end{equation}

The definition of $\mathbf{m}\left(  g_{1},g_{2},\ldots,g_{n}\right)  $
yields
\[
\mathbf{m}\left(  g_{1},g_{2},\ldots,g_{n}\right)  =\dfrac{\left(  g_{1}%
+g_{2}+\cdots+g_{n}\right)  !}{g_{1}!g_{2}!\cdots g_{n}!}.
\]

But the definition of $\mathbf{m}\left(  g_{\sigma\left(  1\right)
},g_{\sigma\left(  2\right)  },\ldots,g_{\sigma\left(  n\right)  }\right)  $
yields%
\begin{align*}
\mathbf{m}\left(  g_{\sigma\left(  1\right)  },g_{\sigma\left(  2\right)
},\ldots,g_{\sigma\left(  n\right)  }\right)   &  =\dfrac{\left(
g_{\sigma\left(  1\right)  }+g_{\sigma\left(  2\right)  }+\cdots
+g_{\sigma\left(  n\right)  }\right)  !}{g_{\sigma\left(  1\right)
}!g_{\sigma\left(  2\right)  }!\cdots g_{\sigma\left(  n\right)  }!}\\
&  =\dfrac{\left(  g_{1}+g_{2}+\cdots+g_{n}\right)  !}{g_{1}!g_{2}!\cdots
g_{n}!}\ \ \ \ \ \ \ \ \ \ \left(  \text{by
(\ref{pf.lem.sol.powerdet.r1.perm-m.short.1}) and
(\ref{pf.lem.sol.powerdet.r1.perm-m.short.2})}\right) \\
&  =\mathbf{m}\left(  g_{1},g_{2},\ldots,g_{n}\right)  .
\end{align*}
This proves Lemma \ref{lem.sol.powerdet.r1.perm-m} \textbf{(b)}.
\end{proof}
\end{vershort}

\begin{verlong}
\begin{proof}
[Proof of Lemma \ref{lem.sol.powerdet.r1.perm-m}.]We have $\sigma\in S_{n}$.
In other words, $\sigma$ is a permutation of the set $\left\{  1,2,\ldots
,n\right\}  $ (since $S_{n}$ is the set of all permutations of the set
$\left\{  1,2,\ldots,n\right\}  $). In other words, $\sigma$ is a bijection
$\left\{  1,2,\ldots,n\right\}  \rightarrow\left\{  1,2,\ldots,n\right\}  $.
Hence, we can substitute $\sigma\left(  i\right)  $ for $i$ in the sum
$\sum_{i\in\left\{  1,2,\ldots,n\right\}  }g_{i}$. We thus obtain%
\[
\sum_{i\in\left\{  1,2,\ldots,n\right\}  }g_{i}=\underbrace{\sum_{i\in\left\{
1,2,\ldots,n\right\}  }}_{=\sum_{i=1}^{n}}g_{\sigma\left(  i\right)  }%
=\sum_{i=1}^{n}g_{\sigma\left(  i\right)  }=g_{\sigma\left(  1\right)
}+g_{\sigma\left(  2\right)  }+\cdots+g_{\sigma\left(  n\right)  }.
\]
Comparing this with%
\[
\underbrace{\sum_{i\in\left\{  1,2,\ldots,n\right\}  }}_{=\sum_{i=1}^{n}}%
g_{i}=\sum_{i=1}^{n}g_{i}=g_{1}+g_{2}+\cdots+g_{n},
\]
we obtain%
\[
g_{\sigma\left(  1\right)  }+g_{\sigma\left(  2\right)  }+\cdots
+g_{\sigma\left(  n\right)  }=g_{1}+g_{2}+\cdots+g_{n}.
\]
This proves Lemma \ref{lem.sol.powerdet.r1.perm-m} \textbf{(a)}.

\textbf{(b)} Recall that $\sigma$ is a bijection $\left\{  1,2,\ldots
,n\right\}  \rightarrow\left\{  1,2,\ldots,n\right\}  $. Hence, we can
substitute $\sigma\left(  i\right)  $ for $i$ in the product $\prod
_{i\in\left\{  1,2,\ldots,n\right\}  }g_{i}!$. We thus obtain%
\[
\prod_{i\in\left\{  1,2,\ldots,n\right\}  }g_{i}!=\underbrace{\prod
_{i\in\left\{  1,2,\ldots,n\right\}  }}_{=\prod_{i=1}^{n}}g_{\sigma\left(
i\right)  }!=\prod_{i=1}^{n}g_{\sigma\left(  i\right)  }!=g_{\sigma\left(
1\right)  }!g_{\sigma\left(  2\right)  }!\cdots g_{\sigma\left(  n\right)
}!.
\]
Comparing this with%
\[
\underbrace{\prod_{i\in\left\{  1,2,\ldots,n\right\}  }}_{=\prod_{i=1}^{n}%
}g_{i}!=\prod_{i=1}^{n}g_{i}!=g_{1}!g_{2}!\cdots g_{n}!,
\]
we obtain%
\[
g_{\sigma\left(  1\right)  }!g_{\sigma\left(  2\right)  }!\cdots
g_{\sigma\left(  n\right)  }!=g_{1}!g_{2}!\cdots g_{n}!.
\]

The definition of $\mathbf{m}\left(  g_{1},g_{2},\ldots,g_{n}\right)  $
yields
\[
\mathbf{m}\left(  g_{1},g_{2},\ldots,g_{n}\right)  =\dfrac{\left(  g_{1}%
+g_{2}+\cdots+g_{n}\right)  !}{g_{1}!g_{2}!\cdots g_{n}!}.
\]

But the definition of $\mathbf{m}\left(  g_{\sigma\left(  1\right)
},g_{\sigma\left(  2\right)  },\ldots,g_{\sigma\left(  n\right)  }\right)  $
yields%
\begin{align*}
\mathbf{m}\left(  g_{\sigma\left(  1\right)  },g_{\sigma\left(  2\right)
},\ldots,g_{\sigma\left(  n\right)  }\right)   &  =\dfrac{\left(
g_{\sigma\left(  1\right)  }+g_{\sigma\left(  2\right)  }+\cdots
+g_{\sigma\left(  n\right)  }\right)  !}{g_{\sigma\left(  1\right)
}!g_{\sigma\left(  2\right)  }!\cdots g_{\sigma\left(  n\right)  }!}\\
&  =\left(  \underbrace{g_{\sigma\left(  1\right)  }+g_{\sigma\left(
2\right)  }+\cdots+g_{\sigma\left(  n\right)  }}_{=g_{1}+g_{2}+\cdots+g_{n}%
}\right)  !/\left(  \underbrace{g_{\sigma\left(  1\right)  }!g_{\sigma\left(
2\right)  }!\cdots g_{\sigma\left(  n\right)  }!}_{=g_{1}!g_{2}!\cdots g_{n}%
!}\right) \\
&  =\left(  g_{1}+g_{2}+\cdots+g_{n}\right)  !/\left(  g_{1}!g_{2}!\cdots
g_{n}!\right) \\
&  =\dfrac{\left(  g_{1}+g_{2}+\cdots+g_{n}\right)  !}{g_{1}!g_{2}!\cdots
g_{n}!}=\mathbf{m}\left(  g_{1},g_{2},\ldots,g_{n}\right)  .
\end{align*}
This proves Lemma \ref{lem.sol.powerdet.r1.perm-m} \textbf{(b)}.
\end{proof}
\end{verlong}

Now, we can attack part \textbf{(c)} of Exercise \ref{exe.powerdet.r1}:

\begin{lemma}
\label{lem.sol.powerdet.r1.c}Let $n\in\mathbb{N}$. Let $a_{1},a_{2}%
,\ldots,a_{n}$ be $n$ elements of $\mathbb{K}$. Let $b_{1},b_{2},\ldots,b_{n}$
be $n$ elements of $\mathbb{K}$.

Let $k\in\mathbb{N}$. Then,%
\begin{align*}
&  \sum_{\sigma\in S_{n}}\left(  -1\right)  ^{\sigma}\left(  \sum_{i=1}%
^{n}a_{i}b_{\sigma\left(  i\right)  }\right)  ^{k}\\
&  =\sum_{\substack{\left(  g_{1},g_{2},\ldots,g_{n}\right)  \in\mathbb{N}%
^{n};\\g_{1}<g_{2}<\cdots<g_{n};\\g_{1}+g_{2}+\cdots+g_{n}=k}}\mathbf{m}%
\left(  g_{1},g_{2},\ldots,g_{n}\right) \\
&  \ \ \ \ \ \ \ \ \ \ \cdot\det\left(  \left(  a_{i}^{g_{j}}\right)  _{1\leq
i\leq n,\ 1\leq j\leq n}\right)  \cdot\det\left(  \left(  b_{i}^{g_{j}%
}\right)  _{1\leq i\leq n,\ 1\leq j\leq n}\right)  .
\end{align*}
Here, we are using the notation introduced in Exercise \ref{exe.multinom2}.
\end{lemma}

\begin{proof}
[Proof of Lemma \ref{lem.sol.powerdet.r1.c}.]We shall use the notations of
Lemma \ref{lem.sol.powerdet.r1.EI}. Recall that%
\begin{align}
\mathbf{E}  &  =\left\{  \left(  k_{1},k_{2},\ldots,k_{n}\right)
\in\mathbb{N}^{n}\ \mid\ \text{the integers }k_{1},k_{2},\ldots,k_{n}\text{
are distinct}\right\} \nonumber\\
&  =\left\{  \left(  g_{1},g_{2},\ldots,g_{n}\right)  \in\mathbb{N}^{n}%
\ \mid\ \text{the integers }g_{1},g_{2},\ldots,g_{n}\text{ are distinct}%
\right\}  \label{pf.lem.sol.powerdet.r1.c.E=}%
\end{align}
(here, we renamed the index $\left(  k_{1},k_{2},\ldots,k_{n}\right)  $ as
$\left(  g_{1},g_{2},\ldots,g_{n}\right)  $) and%
\begin{align}
\mathbf{I}  &  =\left\{  \left(  k_{1},k_{2},\ldots,k_{n}\right)
\in\mathbb{N}^{n}\ \mid\ k_{1}<k_{2}<\cdots<k_{n}\right\} \nonumber\\
&  =\left\{  \left(  g_{1},g_{2},\ldots,g_{n}\right)  \in\mathbb{N}^{n}%
\ \mid\ g_{1}<g_{2}<\cdots<g_{n}\right\}  \label{pf.lem.sol.powerdet.r1.c.I=}%
\end{align}
(here, we renamed the index $\left(  k_{1},k_{2},\ldots,k_{n}\right)  $ as
$\left(  g_{1},g_{2},\ldots,g_{n}\right)  $).

Lemma \ref{lem.sol.powerdet.r1.EI} says that the map%
\begin{align*}
\mathbf{I}\times S_{n}  &  \rightarrow\mathbf{E},\\
\left(  \left(  g_{1},g_{2},\ldots,g_{n}\right)  ,\sigma\right)   &
\mapsto\left(  g_{\sigma\left(  1\right)  },g_{\sigma\left(  2\right)
},\ldots,g_{\sigma\left(  n\right)  }\right)
\end{align*}
is well-defined and is a bijection.

Next, let us make three observations:

\begin{itemize}
\item Every $\left(  k_{1},k_{2},\ldots,k_{n}\right)  \in\mathbb{N}^{n}$
satisfying $\left(  k_{1},k_{2},\ldots,k_{n}\right)  \notin\mathbf{E}$
satisfies%
\begin{equation}
\sum_{\sigma\in S_{n}}\left(  -1\right)  ^{\sigma}\prod_{i=1}^{n}\left(
a_{i}b_{\sigma\left(  i\right)  }\right)  ^{k_{i}}=0
\label{pf.lem.sol.powerdet.r1.c.notEis0}%
\end{equation}
\footnote{\textit{Proof of (\ref{pf.lem.sol.powerdet.r1.c.notEis0}):} Let
$\left(  k_{1},k_{2},\ldots,k_{n}\right)  \in\mathbb{N}^{n}$ be such that
$\left(  k_{1},k_{2},\ldots,k_{n}\right)  \notin\mathbf{E}$.
\par
Let us first show that the integers $k_{1},k_{2},\ldots,k_{n}$ are not
distinct. Indeed, we assume the contrary. Thus, the integers $k_{1}%
,k_{2},\ldots,k_{n}$ are distinct. Hence, $\left(  k_{1},k_{2},\ldots
,k_{n}\right)  $ is an element of $\mathbb{N}^{n}$ having the property that
the integers $k_{1},k_{2},\ldots,k_{n}$ are distinct. In other words, $\left(
k_{1},k_{2},\ldots,k_{n}\right)  $ is an element $\left(  g_{1},g_{2}%
,\ldots,g_{n}\right)  \in\mathbb{N}^{n}$ such that the integers $g_{1}%
,g_{2},\ldots,g_{n}$ are distinct. In other words,%
\[
\left(  k_{1},k_{2},\ldots,k_{n}\right)  \in\left\{  \left(  g_{1}%
,g_{2},\ldots,g_{n}\right)  \in\mathbb{N}^{n}\ \mid\ \text{the integers }%
g_{1},g_{2},\ldots,g_{n}\text{ are distinct}\right\}  .
\]
In light of (\ref{pf.lem.sol.powerdet.r1.c.E=}), this rewrites as $\left(
k_{1},k_{2},\ldots,k_{n}\right)  \in\mathbf{E}$. This contradicts $\left(
k_{1},k_{2},\ldots,k_{n}\right)  \notin\mathbf{E}$. This contradiction proves
that our assumption was wrong. Hence, we have proven that the integers
$k_{1},k_{2},\ldots,k_{n}$ are not distinct. Hence, Lemma
\ref{lem.sol.powerdet.r1.1} \textbf{(b)} yields $\sum_{\sigma\in S_{n}}\left(
-1\right)  ^{\sigma}\prod_{i=1}^{n}\left(  a_{i}b_{\sigma\left(  i\right)
}\right)  ^{k_{i}}=0$. This proves (\ref{pf.lem.sol.powerdet.r1.c.notEis0}).}.

\item Every $\left(  g_{1},g_{2},\ldots,g_{n}\right)  \in\mathbb{N}^{n}$ and
every $\sigma\in S_{n}$ satisfy%
\begin{equation}
\det\left(  \left(  b_{i}^{g_{\sigma\left(  j\right)  }}\right)  _{1\leq i\leq
n,\ 1\leq j\leq n}\right)  =\left(  -1\right)  ^{\sigma}\cdot\det\left(
\left(  b_{i}^{g_{j}}\right)  _{1\leq i\leq n,\ 1\leq j\leq n}\right)
\label{pf.lem.sol.powerdet.r1.c.Eis}%
\end{equation}
\footnote{\textit{Proof of (\ref{pf.lem.sol.powerdet.r1.c.Eis}):} Let $\left(
g_{1},g_{2},\ldots,g_{n}\right)  \in\mathbb{N}^{n}$ and $\sigma\in S_{n}$.
Then, $g_{1},g_{2},\ldots,g_{n}$ are elements of $\mathbb{N}$ (since $\left(
g_{1},g_{2},\ldots,g_{n}\right)  \in\mathbb{N}^{n}$). Hence, $g_{\sigma\left(
1\right)  },g_{\sigma\left(  2\right)  },\ldots,g_{\sigma\left(  n\right)  }$
are elements of $\mathbb{N}$. In other words, $\left(  g_{\sigma\left(
1\right)  },g_{\sigma\left(  2\right)  },\ldots,g_{\sigma\left(  n\right)
}\right)  \in\mathbb{N}^{n}$.
\par
Let $\left[  n\right]  $ denote the set $\left\{  1,2,\ldots,n\right\}  $.
Recall that $\sigma\in S_{n}$. In other words, $\sigma$ is a permutation of
the set $\left\{  1,2,\ldots,n\right\}  $ (since $S_{n}$ is the set of all
permutations of the set $\left\{  1,2,\ldots,n\right\}  $). In other words,
$\sigma$ is a permutation of the set $\left[  n\right]  $ (since $\left\{
1,2,\ldots,n\right\}  =\left[  n\right]  $). In other words, $\sigma$ is a
bijective map $\left[  n\right]  \rightarrow\left[  n\right]  $.
\par
Define an $n\times n$-matrix $A$ by $A=\left(  b_{j}^{g_{i}}\right)  _{1\leq
i\leq n,\ 1\leq j\leq n}$. Define an $n\times n$-matrix $A_{\sigma}$ by
$A_{\sigma}=\left(  b_{j}^{g_{\sigma\left(  i\right)  }}\right)  _{1\leq i\leq
n,\ 1\leq j\leq n}$. Lemma \ref{lem.det.sigma} \textbf{(a)} (applied to
$\sigma$, $A$, $b_{j}^{g_{i}}$ and $A_{\sigma}$ instead of $\kappa$, $B$,
$b_{i,j}$ and $B_{\kappa}$) then yields $\det\left(  A_{\sigma}\right)
=\left(  -1\right)  ^{\sigma}\cdot\det A$.
\par
But $A=\left(  b_{j}^{g_{i}}\right)  _{1\leq i\leq n,\ 1\leq j\leq n}$. Thus,
$A^{T}=\left(  b_{i}^{g_{j}}\right)  _{1\leq i\leq n,\ 1\leq j\leq n}$ (by the
definition of the transpose matrix $A^{T}$). Exercise \ref{exe.ps4.4} now
yields $\det\left(  A^{T}\right)  =\det A$. Hence,%
\[
\det A=\det\left(  \underbrace{A^{T}}_{=\left(  b_{i}^{g_{j}}\right)  _{1\leq
i\leq n,\ 1\leq j\leq n}}\right)  =\det\left(  \left(  b_{i}^{g_{j}}\right)
_{1\leq i\leq n,\ 1\leq j\leq n}\right)  .
\]
\par
On the other hand, $A_{\sigma}=\left(  b_{j}^{g_{\sigma\left(  i\right)  }%
}\right)  _{1\leq i\leq n,\ 1\leq j\leq n}$. Hence, $\left(  A_{\sigma
}\right)  ^{T}=\left(  b_{i}^{g_{\sigma\left(  j\right)  }}\right)  _{1\leq
i\leq n,\ 1\leq j\leq n}$ (by the definition of the transpose matrix $\left(
A_{\sigma}\right)  ^{T}$). Exercise \ref{exe.ps4.4} (applied to $A_{\sigma}$
instead of $A$) now yields $\det\left(  \left(  A_{\sigma}\right)
^{T}\right)  =\det\left(  A_{\sigma}\right)  $. Since $\left(  A_{\sigma
}\right)  ^{T}=\left(  b_{i}^{g_{\sigma\left(  j\right)  }}\right)  _{1\leq
i\leq n,\ 1\leq j\leq n}$, this rewrites as%
\begin{align*}
\det\left(  \left(  b_{i}^{g_{\sigma\left(  j\right)  }}\right)  _{1\leq i\leq
n,\ 1\leq j\leq n}\right)   &  =\det\left(  A_{\sigma}\right)  =\left(
-1\right)  ^{\sigma}\cdot\underbrace{\det A}_{=\det\left(  \left(
b_{i}^{g_{j}}\right)  _{1\leq i\leq n,\ 1\leq j\leq n}\right)  }\\
&  =\left(  -1\right)  ^{\sigma}\cdot\det\left(  \left(  b_{i}^{g_{j}}\right)
_{1\leq i\leq n,\ 1\leq j\leq n}\right)  .
\end{align*}
This proves (\ref{pf.lem.sol.powerdet.r1.c.Eis}).}.

\item Every $\left(  g_{1},g_{2},\ldots,g_{n}\right)  \in\mathbb{N}^{n}$
satisfies%
\begin{equation}
\det\left(  \left(  a_{i}^{g_{j}}\right)  _{1\leq i\leq n,\ 1\leq j\leq
n}\right)  =\sum_{\sigma\in S_{n}}\left(  \prod_{i=1}^{n}a_{i}^{g_{\sigma
\left(  i\right)  }}\right)  \left(  -1\right)  ^{\sigma}
\label{pf.lem.sol.powerdet.r1.c.adet}%
\end{equation}
\footnote{\textit{Proof of (\ref{pf.lem.sol.powerdet.r1.c.adet}):} Let
$\left(  g_{1},g_{2},\ldots,g_{n}\right)  \in\mathbb{N}^{n}$. The equality
(\ref{eq.det.eq.2}) (applied to $\left(  a_{i}^{g_{j}}\right)  _{1\leq i\leq
n,\ 1\leq j\leq n}$ and $a_{i}^{g_{j}}$ instead of $A$ and $a_{i,j}$) yields%
\[
\det\left(  \left(  a_{i}^{g_{j}}\right)  _{1\leq i\leq n,\ 1\leq j\leq
n}\right)  =\sum_{\sigma\in S_{n}}\left(  -1\right)  ^{\sigma}\prod_{i=1}%
^{n}a_{i}^{g_{\sigma\left(  i\right)  }}=\sum_{\sigma\in S_{n}}\left(
\prod_{i=1}^{n}a_{i}^{g_{\sigma\left(  i\right)  }}\right)  \left(  -1\right)
^{\sigma}.
\]
This proves (\ref{pf.lem.sol.powerdet.r1.c.adet}).}.
\end{itemize}

Now, each $\sigma\in S_{n}$ satisfies%
\begin{equation}
\left(  \sum_{i=1}^{n}a_{i}b_{\sigma\left(  i\right)  }\right)  ^{k}%
=\sum_{\substack{\left(  k_{1},k_{2},\ldots,k_{n}\right)  \in\mathbb{N}%
^{n};\\k_{1}+k_{2}+\cdots+k_{n}=k}}\mathbf{m}\left(  k_{1},k_{2},\ldots
,k_{n}\right)  \prod_{i=1}^{n}\left(  a_{i}b_{\sigma\left(  i\right)
}\right)  ^{k_{i}} \label{pf.lem.sol.powerdet.r1.c.mulnom}%
\end{equation}
\footnote{\textit{Proof of (\ref{pf.lem.sol.powerdet.r1.c.mulnom}):} Let
$\sigma\in S_{n}$. Then,%
\begin{align*}
\left(  \underbrace{\sum_{i=1}^{n}a_{i}b_{\sigma\left(  i\right)  }}%
_{=a_{1}b_{\sigma\left(  1\right)  }+a_{2}b_{\sigma\left(  2\right)  }%
+\cdots+a_{n}b_{\sigma\left(  n\right)  }}\right)  ^{k}  &  =\left(
a_{1}b_{\sigma\left(  1\right)  }+a_{2}b_{\sigma\left(  2\right)  }%
+\cdots+a_{n}b_{\sigma\left(  n\right)  }\right)  ^{k}\\
&  =\sum_{\substack{\left(  k_{1},k_{2},\ldots,k_{n}\right)  \in\mathbb{N}%
^{n};\\k_{1}+k_{2}+\cdots+k_{n}=k}}\mathbf{m}\left(  k_{1},k_{2},\ldots
,k_{n}\right)  \prod_{i=1}^{n}\left(  a_{i}b_{\sigma\left(  i\right)
}\right)  ^{k_{i}}%
\end{align*}
(by Exercise \ref{exe.multinom2} (applied to $n$, $k$ and $a_{i}%
b_{\sigma\left(  i\right)  }$ instead of $m$, $n$ and $a_{i}$)). This proves
(\ref{pf.lem.sol.powerdet.r1.c.mulnom}).}.

Hence,%
\begin{align*}
&  \sum_{\sigma\in S_{n}}\left(  -1\right)  ^{\sigma}\underbrace{\left(
\sum_{i=1}^{n}a_{i}b_{\sigma\left(  i\right)  }\right)  ^{k}}_{\substack{=\sum
_{\substack{\left(  k_{1},k_{2},\ldots,k_{n}\right)  \in\mathbb{N}^{n}%
;\\k_{1}+k_{2}+\cdots+k_{n}=k}}\mathbf{m}\left(  k_{1},k_{2},\ldots
,k_{n}\right)  \prod_{i=1}^{n}\left(  a_{i}b_{\sigma\left(  i\right)
}\right)  ^{k_{i}}\\\text{(by (\ref{pf.lem.sol.powerdet.r1.c.mulnom}))}}}\\
&  =\sum_{\sigma\in S_{n}}\left(  -1\right)  ^{\sigma}\sum_{\substack{\left(
k_{1},k_{2},\ldots,k_{n}\right)  \in\mathbb{N}^{n};\\k_{1}+k_{2}+\cdots
+k_{n}=k}}\mathbf{m}\left(  k_{1},k_{2},\ldots,k_{n}\right)  \prod_{i=1}%
^{n}\left(  a_{i}b_{\sigma\left(  i\right)  }\right)  ^{k_{i}}\\
&  =\underbrace{\sum_{\sigma\in S_{n}}\sum_{\substack{\left(  k_{1}%
,k_{2},\ldots,k_{n}\right)  \in\mathbb{N}^{n};\\k_{1}+k_{2}+\cdots+k_{n}=k}%
}}_{=\sum_{\substack{\left(  k_{1},k_{2},\ldots,k_{n}\right)  \in
\mathbb{N}^{n};\\k_{1}+k_{2}+\cdots+k_{n}=k}}\sum_{\sigma\in S_{n}}}\left(
-1\right)  ^{\sigma}\mathbf{m}\left(  k_{1},k_{2},\ldots,k_{n}\right)
\prod_{i=1}^{n}\left(  a_{i}b_{\sigma\left(  i\right)  }\right)  ^{k_{i}}\\
&  =\sum_{\substack{\left(  k_{1},k_{2},\ldots,k_{n}\right)  \in\mathbb{N}%
^{n};\\k_{1}+k_{2}+\cdots+k_{n}=k}}\sum_{\sigma\in S_{n}}\left(  -1\right)
^{\sigma}\mathbf{m}\left(  k_{1},k_{2},\ldots,k_{n}\right)  \prod_{i=1}%
^{n}\left(  a_{i}b_{\sigma\left(  i\right)  }\right)  ^{k_{i}}\\
&  =\sum_{\substack{\left(  k_{1},k_{2},\ldots,k_{n}\right)  \in\mathbb{N}%
^{n};\\k_{1}+k_{2}+\cdots+k_{n}=k}}\mathbf{m}\left(  k_{1},k_{2},\ldots
,k_{n}\right)  \sum_{\sigma\in S_{n}}\left(  -1\right)  ^{\sigma}\prod
_{i=1}^{n}\left(  a_{i}b_{\sigma\left(  i\right)  }\right)  ^{k_{i}}\\
&  =\underbrace{\sum_{\substack{\left(  k_{1},k_{2},\ldots,k_{n}\right)
\in\mathbb{N}^{n};\\k_{1}+k_{2}+\cdots+k_{n}=k;\\\left(  k_{1},k_{2}%
,\ldots,k_{n}\right)  \in\mathbf{E}}}}_{\substack{=\sum_{\substack{\left(
k_{1},k_{2},\ldots,k_{n}\right)  \in\mathbb{N}^{n};\\\left(  k_{1}%
,k_{2},\ldots,k_{n}\right)  \in\mathbf{E};\\k_{1}+k_{2}+\cdots+k_{n}%
=k}}\\=\sum_{\substack{\left(  k_{1},k_{2},\ldots,k_{n}\right)  \in
\mathbf{E};\\k_{1}+k_{2}+\cdots+k_{n}=k}}\\\text{(since }\mathbf{E}%
\subseteq\mathbb{N}^{n}\text{)}}}\mathbf{m}\left(  k_{1},k_{2},\ldots
,k_{n}\right)  \underbrace{\sum_{\sigma\in S_{n}}\left(  -1\right)  ^{\sigma
}\prod_{i=1}^{n}\left(  a_{i}b_{\sigma\left(  i\right)  }\right)  ^{k_{i}}%
}_{\substack{=\left(  \prod_{i=1}^{n}a_{i}^{k_{i}}\right)  \det\left(  \left(
b_{i}^{k_{j}}\right)  _{1\leq i\leq n,\ 1\leq j\leq n}\right)  \\\text{(by
Lemma \ref{lem.sol.powerdet.r1.1} \textbf{(a)})}}}\\
&  \ \ \ \ \ \ \ \ \ \ +\sum_{\substack{\left(  k_{1},k_{2},\ldots
,k_{n}\right)  \in\mathbb{N}^{n};\\k_{1}+k_{2}+\cdots+k_{n}=k;\\\left(
k_{1},k_{2},\ldots,k_{n}\right)  \notin\mathbf{E}}}\mathbf{m}\left(
k_{1},k_{2},\ldots,k_{n}\right)  \underbrace{\sum_{\sigma\in S_{n}}\left(
-1\right)  ^{\sigma}\prod_{i=1}^{n}\left(  a_{i}b_{\sigma\left(  i\right)
}\right)  ^{k_{i}}}_{\substack{=0\\\text{(by
(\ref{pf.lem.sol.powerdet.r1.c.notEis0}))}}}\\
&  \ \ \ \ \ \ \ \ \ \ \left(
\begin{array}
[c]{c}%
\text{since each }\left(  k_{1},k_{2},\ldots,k_{n}\right)  \in\mathbb{N}%
^{n}\text{ satisfies either }\left(  k_{1},k_{2},\ldots,k_{n}\right)
\in\mathbf{E}\\
\text{or }\left(  k_{1},k_{2},\ldots,k_{n}\right)  \notin\mathbf{E}\text{ (but
not both)}%
\end{array}
\right)
\end{align*}%
\begin{align}
&  =\sum_{\substack{\left(  k_{1},k_{2},\ldots,k_{n}\right)  \in
\mathbf{E};\\k_{1}+k_{2}+\cdots+k_{n}=k}}\mathbf{m}\left(  k_{1},k_{2}%
,\ldots,k_{n}\right)  \left(  \prod_{i=1}^{n}a_{i}^{k_{i}}\right)  \det\left(
\left(  b_{i}^{k_{j}}\right)  _{1\leq i\leq n,\ 1\leq j\leq n}\right)
\nonumber\\
&  \ \ \ \ \ \ \ \ \ \ +\underbrace{\sum_{\substack{\left(  k_{1},k_{2}%
,\ldots,k_{n}\right)  \in\mathbb{N}^{n};\\k_{1}+k_{2}+\cdots+k_{n}=k;\\\left(
k_{1},k_{2},\ldots,k_{n}\right)  \notin\mathbf{E}}}\mathbf{m}\left(
k_{1},k_{2},\ldots,k_{n}\right)  0}_{=0}\nonumber\\
&  =\sum_{\substack{\left(  k_{1},k_{2},\ldots,k_{n}\right)  \in
\mathbf{E};\\k_{1}+k_{2}+\cdots+k_{n}=k}}\mathbf{m}\left(  k_{1},k_{2}%
,\ldots,k_{n}\right)  \left(  \prod_{i=1}^{n}a_{i}^{k_{i}}\right)  \det\left(
\left(  b_{i}^{k_{j}}\right)  _{1\leq i\leq n,\ 1\leq j\leq n}\right)
\nonumber\\
&  =\sum_{\substack{\left(  \left(  g_{1},g_{2},\ldots,g_{n}\right)
,\sigma\right)  \in\mathbf{I}\times S_{n};\\g_{\sigma\left(  1\right)
}+g_{\sigma\left(  2\right)  }+\cdots+g_{\sigma\left(  n\right)  }%
=k}}\mathbf{m}\left(  g_{\sigma\left(  1\right)  },g_{\sigma\left(  2\right)
},\ldots,g_{\sigma\left(  n\right)  }\right)  \left(  \prod_{i=1}^{n}%
a_{i}^{g_{\sigma\left(  i\right)  }}\right)  \det\left(  \left(
b_{i}^{g_{\sigma\left(  j\right)  }}\right)  _{1\leq i\leq n,\ 1\leq j\leq
n}\right)  \label{pf.lem.sol.powerdet.r1.c.3}%
\end{align}
(here, we have substituted $\left(  g_{\sigma\left(  1\right)  }%
,g_{\sigma\left(  2\right)  },\ldots,g_{\sigma\left(  n\right)  }\right)  $
for $\left(  k_{1},k_{2},\ldots,k_{n}\right)  $ in the sum, because the map%
\begin{align*}
\mathbf{I}\times S_{n}  &  \rightarrow\mathbf{E},\\
\left(  \left(  g_{1},g_{2},\ldots,g_{n}\right)  ,\sigma\right)   &
\mapsto\left(  g_{\sigma\left(  1\right)  },g_{\sigma\left(  2\right)
},\ldots,g_{\sigma\left(  n\right)  }\right)
\end{align*}
is a bijection).

But we have the following equality of summation signs:%
\begin{align*}
\sum_{\substack{\left(  \left(  g_{1},g_{2},\ldots,g_{n}\right)
,\sigma\right)  \in\mathbf{I}\times S_{n};\\g_{\sigma\left(  1\right)
}+g_{\sigma\left(  2\right)  }+\cdots+g_{\sigma\left(  n\right)  }=k}}  &
=\underbrace{\sum_{\left(  g_{1},g_{2},\ldots,g_{n}\right)  \in\mathbf{I}}%
}_{\substack{=\sum_{\substack{\left(  g_{1},g_{2},\ldots,g_{n}\right)
\in\mathbb{N}^{n};\\g_{1}<g_{2}<\cdots<g_{n}}}\\\text{(since }\mathbf{I}%
=\left\{  \left(  g_{1},g_{2},\ldots,g_{n}\right)  \in\mathbb{N}^{n}%
\ \mid\ g_{1}<g_{2}<\cdots<g_{n}\right\}  \text{)}}}\underbrace{\sum
_{\substack{\sigma\in S_{n};\\g_{\sigma\left(  1\right)  }+g_{\sigma\left(
2\right)  }+\cdots+g_{\sigma\left(  n\right)  }=k}}}_{\substack{=\sum
_{\substack{\sigma\in S_{n};\\g_{1}+g_{2}+\cdots+g_{n}=k}}\\\text{(since each
}\sigma\in S_{n}\text{ satisfies}\\g_{\sigma\left(  1\right)  }+g_{\sigma
\left(  2\right)  }+\cdots+g_{\sigma\left(  n\right)  }=g_{1}+g_{2}%
+\cdots+g_{n}\\\text{(by Lemma \ref{lem.sol.powerdet.r1.perm-m} \textbf{(a)}%
))}}}\\
&  =\sum_{\substack{\left(  g_{1},g_{2},\ldots,g_{n}\right)  \in\mathbb{N}%
^{n};\\g_{1}<g_{2}<\cdots<g_{n}}}\sum_{\substack{\sigma\in S_{n};\\g_{1}%
+g_{2}+\cdots+g_{n}=k}}=\sum_{\substack{\left(  g_{1},g_{2},\ldots
,g_{n}\right)  \in\mathbb{N}^{n};\\g_{1}<g_{2}<\cdots<g_{n};\\g_{1}%
+g_{2}+\cdots+g_{n}=k}}\sum_{\sigma\in S_{n}}.
\end{align*}
Hence, (\ref{pf.lem.sol.powerdet.r1.c.3}) becomes%
\begin{align*}
&  \sum_{\sigma\in S_{n}}\left(  -1\right)  ^{\sigma}\left(  \sum_{i=1}%
^{n}a_{i}b_{\sigma\left(  i\right)  }\right)  ^{k}\\
&  =\underbrace{\sum_{\substack{\left(  \left(  g_{1},g_{2},\ldots
,g_{n}\right)  ,\sigma\right)  \in\mathbf{I}\times S_{n};\\g_{\sigma\left(
1\right)  }+g_{\sigma\left(  2\right)  }+\cdots+g_{\sigma\left(  n\right)
}=k}}}_{=\sum_{\substack{\left(  g_{1},g_{2},\ldots,g_{n}\right)
\in\mathbb{N}^{n};\\g_{1}<g_{2}<\cdots<g_{n};\\g_{1}+g_{2}+\cdots+g_{n}%
=k}}\sum_{\sigma\in S_{n}}}\mathbf{m}\left(  g_{\sigma\left(  1\right)
},g_{\sigma\left(  2\right)  },\ldots,g_{\sigma\left(  n\right)  }\right)
\left(  \prod_{i=1}^{n}a_{i}^{g_{\sigma\left(  i\right)  }}\right)
\det\left(  \left(  b_{i}^{g_{\sigma\left(  j\right)  }}\right)  _{1\leq i\leq
n,\ 1\leq j\leq n}\right) \\
&  =\sum_{\substack{\left(  g_{1},g_{2},\ldots,g_{n}\right)  \in\mathbb{N}%
^{n};\\g_{1}<g_{2}<\cdots<g_{n};\\g_{1}+g_{2}+\cdots+g_{n}=k}}\sum_{\sigma\in
S_{n}}\underbrace{\mathbf{m}\left(  g_{\sigma\left(  1\right)  }%
,g_{\sigma\left(  2\right)  },\ldots,g_{\sigma\left(  n\right)  }\right)
}_{\substack{=\mathbf{m}\left(  g_{1},g_{2},\ldots,g_{n}\right)  \\\text{(by
Lemma \ref{lem.sol.powerdet.r1.perm-m} \textbf{(b)})}}}\left(  \prod_{i=1}%
^{n}a_{i}^{g_{\sigma\left(  i\right)  }}\right)  \underbrace{\det\left(
\left(  b_{i}^{g_{\sigma\left(  j\right)  }}\right)  _{1\leq i\leq n,\ 1\leq
j\leq n}\right)  }_{\substack{=\left(  -1\right)  ^{\sigma}\cdot\det\left(
\left(  b_{i}^{g_{j}}\right)  _{1\leq i\leq n,\ 1\leq j\leq n}\right)
\\\text{(by (\ref{pf.lem.sol.powerdet.r1.c.Eis}))}}}\\
&  =\sum_{\substack{\left(  g_{1},g_{2},\ldots,g_{n}\right)  \in\mathbb{N}%
^{n};\\g_{1}<g_{2}<\cdots<g_{n};\\g_{1}+g_{2}+\cdots+g_{n}=k}}\underbrace{\sum
_{\sigma\in S_{n}}\mathbf{m}\left(  g_{1},g_{2},\ldots,g_{n}\right)  \left(
\prod_{i=1}^{n}a_{i}^{g_{\sigma\left(  i\right)  }}\right)  \left(  -1\right)
^{\sigma}\cdot\det\left(  \left(  b_{i}^{g_{j}}\right)  _{1\leq i\leq
n,\ 1\leq j\leq n}\right)  }_{=\mathbf{m}\left(  g_{1},g_{2},\ldots
,g_{n}\right)  \cdot\left(  \sum_{\sigma\in S_{n}}\left(  \prod_{i=1}^{n}%
a_{i}^{g_{\sigma\left(  i\right)  }}\right)  \left(  -1\right)  ^{\sigma
}\right)  \det\left(  \left(  b_{i}^{g_{j}}\right)  _{1\leq i\leq n,\ 1\leq
j\leq n}\right)  }\\
&  =\sum_{\substack{\left(  g_{1},g_{2},\ldots,g_{n}\right)  \in\mathbb{N}%
^{n};\\g_{1}<g_{2}<\cdots<g_{n};\\g_{1}+g_{2}+\cdots+g_{n}=k}}\mathbf{m}%
\left(  g_{1},g_{2},\ldots,g_{n}\right)  \cdot\underbrace{\left(  \sum
_{\sigma\in S_{n}}\left(  \prod_{i=1}^{n}a_{i}^{g_{\sigma\left(  i\right)  }%
}\right)  \left(  -1\right)  ^{\sigma}\right)  }_{\substack{=\det\left(
\left(  a_{i}^{g_{j}}\right)  _{1\leq i\leq n,\ 1\leq j\leq n}\right)
\\\text{(by (\ref{pf.lem.sol.powerdet.r1.c.adet}))}}}\det\left(  \left(
b_{i}^{g_{j}}\right)  _{1\leq i\leq n,\ 1\leq j\leq n}\right) \\
&  =\sum_{\substack{\left(  g_{1},g_{2},\ldots,g_{n}\right)  \in\mathbb{N}%
^{n};\\g_{1}<g_{2}<\cdots<g_{n};\\g_{1}+g_{2}+\cdots+g_{n}=k}}\mathbf{m}%
\left(  g_{1},g_{2},\ldots,g_{n}\right)  \cdot\det\left(  \left(  a_{i}%
^{g_{j}}\right)  _{1\leq i\leq n,\ 1\leq j\leq n}\right)  \cdot\det\left(
\left(  b_{i}^{g_{j}}\right)  _{1\leq i\leq n,\ 1\leq j\leq n}\right)  .
\end{align*}
This proves Lemma \ref{lem.sol.powerdet.r1.c}.
\end{proof}

Next, we can use Lemma \ref{lem.sol.powerdet.r1.c} to solve part \textbf{(a)}
of Exercise \ref{exe.powerdet.r1}:

\begin{lemma}
\label{lem.sol.powerdet.r1.a}Let $n\in\mathbb{N}$. Let $a_{1},a_{2}%
,\ldots,a_{n}$ be $n$ elements of $\mathbb{K}$. Let $b_{1},b_{2},\ldots,b_{n}$
be $n$ elements of $\mathbb{K}$. Let $m=\dbinom{n}{2}$.

Then,
\[
\sum_{\sigma\in S_{n}}\left(  -1\right)  ^{\sigma}\left(  \sum_{i=1}^{n}%
a_{i}b_{\sigma\left(  i\right)  }\right)  ^{k}=0
\]
for each $k\in\left\{  0,1,\ldots,m-1\right\}  $.
\end{lemma}

\begin{proof}
[Proof of Lemma \ref{lem.sol.powerdet.r1.a}.]Let $k\in\left\{  0,1,\ldots
,m-1\right\}  $. Then, $k\in\mathbb{N}$ and $k\leq m-1$.

There exists no $\left(  g_{1},g_{2},\ldots,g_{n}\right)  \in\mathbb{N}^{n}$
satisfying $g_{1}<g_{2}<\cdots<g_{n}$ and $g_{1}+g_{2}+\cdots+g_{n}%
=k$\ \ \ \ \footnote{\textit{Proof.} Let $\left(  g_{1},g_{2},\ldots
,g_{n}\right)  \in\mathbb{N}^{n}$ be such that $g_{1}<g_{2}<\cdots<g_{n}$ and
$g_{1}+g_{2}+\cdots+g_{n}=k$. We shall derive a contradiction.
\par
Lemma \ref{lem.sol.powerdet.r1.0} \textbf{(a)} (applied to $k_{i}=g_{i}$)
yields $g_{1}+g_{2}+\cdots+g_{n}\geq\dbinom{n}{2}$. This contradicts
$g_{1}+g_{2}+\cdots+g_{n}=k\leq m-1<m=\dbinom{n}{2}$.
\par
Now, let us forget that we fixed $\left(  g_{1},g_{2},\ldots,g_{n}\right)  $.
We thus have derived a contradiction for each $\left(  g_{1},g_{2}%
,\ldots,g_{n}\right)  \in\mathbb{N}^{n}$ satisfying $g_{1}<g_{2}<\cdots<g_{n}$
and $g_{1}+g_{2}+\cdots+g_{n}=k$. Hence, there exists no $\left(  g_{1}%
,g_{2},\ldots,g_{n}\right)  \in\mathbb{N}^{n}$ satisfying $g_{1}<g_{2}%
<\cdots<g_{n}$ and $g_{1}+g_{2}+\cdots+g_{n}=k$.}. Thus, the sum
\[
\sum_{\substack{\left(  g_{1},g_{2},\ldots,g_{n}\right)  \in\mathbb{N}%
^{n};\\g_{1}<g_{2}<\cdots<g_{n};\\g_{1}+g_{2}+\cdots+g_{n}=k}}\mathbf{m}%
\left(  g_{1},g_{2},\ldots,g_{n}\right)  \cdot\det\left(  \left(  a_{i}%
^{g_{j}}\right)  _{1\leq i\leq n,\ 1\leq j\leq n}\right)  \cdot\det\left(
\left(  b_{i}^{g_{j}}\right)  _{1\leq i\leq n,\ 1\leq j\leq n}\right)
\]
is an empty sum. Hence,%
\begin{align*}
&  \sum_{\substack{\left(  g_{1},g_{2},\ldots,g_{n}\right)  \in\mathbb{N}%
^{n};\\g_{1}<g_{2}<\cdots<g_{n};\\g_{1}+g_{2}+\cdots+g_{n}=k}}\mathbf{m}%
\left(  g_{1},g_{2},\ldots,g_{n}\right)  \cdot\det\left(  \left(  a_{i}%
^{g_{j}}\right)  _{1\leq i\leq n,\ 1\leq j\leq n}\right)  \cdot\det\left(
\left(  b_{i}^{g_{j}}\right)  _{1\leq i\leq n,\ 1\leq j\leq n}\right) \\
&  =\left(  \text{empty sum}\right)  =0.
\end{align*}

Now, Lemma \ref{lem.sol.powerdet.r1.c} yields%
\begin{align*}
&  \sum_{\sigma\in S_{n}}\left(  -1\right)  ^{\sigma}\left(  \sum_{i=1}%
^{n}a_{i}b_{\sigma\left(  i\right)  }\right)  ^{k}\\
&  =\sum_{\substack{\left(  g_{1},g_{2},\ldots,g_{n}\right)  \in\mathbb{N}%
^{n};\\g_{1}<g_{2}<\cdots<g_{n};\\g_{1}+g_{2}+\cdots+g_{n}=k}}\mathbf{m}%
\left(  g_{1},g_{2},\ldots,g_{n}\right)  \cdot\det\left(  \left(  a_{i}%
^{g_{j}}\right)  _{1\leq i\leq n,\ 1\leq j\leq n}\right)  \cdot\det\left(
\left(  b_{i}^{g_{j}}\right)  _{1\leq i\leq n,\ 1\leq j\leq n}\right) \\
&  =0.
\end{align*}
This proves Lemma \ref{lem.sol.powerdet.r1.a}.
\end{proof}

Finally, we are ready for part \textbf{(b)} of Exercise \ref{exe.powerdet.r1}:

\begin{lemma}
\label{lem.sol.powerdet.r1.b}Let $n\in\mathbb{N}$. Let $a_{1},a_{2}%
,\ldots,a_{n}$ be $n$ elements of $\mathbb{K}$. Let $b_{1},b_{2},\ldots,b_{n}$
be $n$ elements of $\mathbb{K}$.

Let $m=\dbinom{n}{2}$. Then,%
\[
\sum_{\sigma\in S_{n}}\left(  -1\right)  ^{\sigma}\left(  \sum_{i=1}^{n}%
a_{i}b_{\sigma\left(  i\right)  }\right)  ^{m}=\mathbf{m}\left(
0,1,\ldots,n-1\right)  \cdot\prod_{1\leq i<j\leq n}\left(  \left(  a_{i}%
-a_{j}\right)  \left(  b_{i}-b_{j}\right)  \right)  .
\]
Here, we are using the notation introduced in Exercise \ref{exe.multinom2}.
\end{lemma}

\begin{proof}
[Proof of Lemma \ref{lem.sol.powerdet.r1.b}.]We have $m=\dbinom{n}{2}%
=\sum_{r=0}^{n-1}r$ (by Lemma \ref{lem.sol.vander-hook.gauss}), thus
$m=\sum_{r=0}^{n-1}r\in\mathbb{N}$. Hence, Lemma \ref{lem.sol.powerdet.r1.c}
(applied to $k=m$) yields%
\begin{align*}
&  \sum_{\sigma\in S_{n}}\left(  -1\right)  ^{\sigma}\left(  \sum_{i=1}%
^{n}a_{i}b_{\sigma\left(  i\right)  }\right)  ^{m}\\
&  =\sum_{\substack{\left(  g_{1},g_{2},\ldots,g_{n}\right)  \in\mathbb{N}%
^{n};\\g_{1}<g_{2}<\cdots<g_{n};\\g_{1}+g_{2}+\cdots+g_{n}=m}}\mathbf{m}%
\left(  g_{1},g_{2},\ldots,g_{n}\right)  \cdot\det\left(  \left(  a_{i}%
^{g_{j}}\right)  _{1\leq i\leq n,\ 1\leq j\leq n}\right)  \cdot\det\left(
\left(  b_{i}^{g_{j}}\right)  _{1\leq i\leq n,\ 1\leq j\leq n}\right)  .
\end{align*}

But we have the following two observations:

\begin{itemize}
\item The $n$-tuple $\left(  0,1,\ldots,n-1\right)  $ is an $n$-tuple $\left(
g_{1},g_{2},\ldots,g_{n}\right)  \in\mathbb{N}^{n}$ satisfying $g_{1}%
<g_{2}<\cdots<g_{n}$ and $g_{1}+g_{2}+\cdots+g_{n}=m$%
\ \ \ \ \footnote{\textit{Proof.} The $n$-tuple $\left(  0,1,\ldots
,n-1\right)  $ belongs to $\mathbb{N}^{n}$ and satisfies $0<1<\cdots<n-1$ and
$0+1+\cdots+\left(  n-1\right)  =m$ (since $0+1+\cdots+\left(  n-1\right)
=\sum_{r=0}^{n-1}r=m$). In other words, the $n$-tuple $\left(  0,1,\ldots
,n-1\right)  $ is an $n$-tuple $\left(  g_{1},g_{2},\ldots,g_{n}\right)
\in\mathbb{N}^{n}$ satisfying $g_{1}<g_{2}<\cdots<g_{n}$ and $g_{1}%
+g_{2}+\cdots+g_{n}=m$.}.

\item Any $n$-tuple $\left(  g_{1},g_{2},\ldots,g_{n}\right)  \in
\mathbb{N}^{n}$ satisfying $g_{1}<g_{2}<\cdots<g_{n}$ and $g_{1}+g_{2}%
+\cdots+g_{n}=m$\ must be equal to $\left(  0,1,\ldots,n-1\right)
$\ \ \ \ \footnote{\textit{Proof.} Let $\left(  g_{1},g_{2},\ldots
,g_{n}\right)  \in\mathbb{N}^{n}$ be any $n$-tuple satisfying $g_{1}%
<g_{2}<\cdots<g_{n}$ and $g_{1}+g_{2}+\cdots+g_{n}=m$. We shall show that
$\left(  g_{1},g_{2},\ldots,g_{n}\right)  $ must be equal to $\left(
0,1,\ldots,n-1\right)  $.
\par
We have $g_{1}+g_{2}+\cdots+g_{n}=m\leq m=\dbinom{n}{2}$. Hence, Lemma
\ref{lem.sol.powerdet.r1.0} \textbf{(b)} (applied to $g_{i}$ instead of
$k_{i}$) shows that $\left(  g_{1},g_{2},\ldots,g_{n}\right)  =\left(
0,1,\ldots,n-1\right)  $. In other words, $\left(  g_{1},g_{2},\ldots
,g_{n}\right)  $ must be equal to $\left(  0,1,\ldots,n-1\right)  $. Qed.}.
\end{itemize}

Combining these two observations, we conclude that there exists exactly one
$n$-tuple $\left(  g_{1},g_{2},\ldots,g_{n}\right)  \in\mathbb{N}^{n}$
satisfying $g_{1}<g_{2}<\cdots<g_{n}$ and $g_{1}+g_{2}+\cdots+g_{n}=m$:
namely, the $n$-tuple $\left(  0,1,\ldots,n-1\right)  $. Hence, the sum%
\[
\sum_{\substack{\left(  g_{1},g_{2},\ldots,g_{n}\right)  \in\mathbb{N}%
^{n};\\g_{1}<g_{2}<\cdots<g_{n};\\g_{1}+g_{2}+\cdots+g_{n}=m}}\mathbf{m}%
\left(  g_{1},g_{2},\ldots,g_{n}\right)  \cdot\det\left(  \left(  a_{i}%
^{g_{j}}\right)  _{1\leq i\leq n,\ 1\leq j\leq n}\right)  \cdot\det\left(
\left(  b_{i}^{g_{j}}\right)  _{1\leq i\leq n,\ 1\leq j\leq n}\right)
\]
has only one addend: namely, the addend for $\left(  g_{1},g_{2},\ldots
,g_{n}\right)  =\left(  0,1,\ldots,n-1\right)  $. Therefore, this sum can be
simplified as follows:%
\begin{align*}
&  \sum_{\substack{\left(  g_{1},g_{2},\ldots,g_{n}\right)  \in\mathbb{N}%
^{n};\\g_{1}<g_{2}<\cdots<g_{n};\\g_{1}+g_{2}+\cdots+g_{n}=m}}\mathbf{m}%
\left(  g_{1},g_{2},\ldots,g_{n}\right)  \cdot\det\left(  \left(  a_{i}%
^{g_{j}}\right)  _{1\leq i\leq n,\ 1\leq j\leq n}\right)  \cdot\det\left(
\left(  b_{i}^{g_{j}}\right)  _{1\leq i\leq n,\ 1\leq j\leq n}\right) \\
&  =\mathbf{m}\left(  0,1,\ldots,n-1\right)  \cdot\underbrace{\det\left(
\left(  a_{i}^{j-1}\right)  _{1\leq i\leq n,\ 1\leq j\leq n}\right)
}_{\substack{=\prod_{1\leq j<i\leq n}\left(  a_{i}-a_{j}\right)  \\\text{(by
Theorem \ref{thm.vander-det} \textbf{(c)}}\\\text{(applied to }x_{k}%
=a_{k}\text{))}}}\cdot\underbrace{\det\left(  \left(  b_{i}^{j-1}\right)
_{1\leq i\leq n,\ 1\leq j\leq n}\right)  }_{\substack{=\prod_{1\leq j<i\leq
n}\left(  b_{i}-b_{j}\right)  \\\text{(by Theorem \ref{thm.vander-det}
\textbf{(c)}}\\\text{(applied to }x_{k}=b_{k}\text{))}}}\\
&  =\mathbf{m}\left(  0,1,\ldots,n-1\right)  \cdot\underbrace{\left(
\prod_{1\leq j<i\leq n}\left(  a_{i}-a_{j}\right)  \right)  \cdot\left(
\prod_{1\leq j<i\leq n}\left(  b_{i}-b_{j}\right)  \right)  }%
_{\substack{=\prod_{1\leq j<i\leq n}\left(  \left(  a_{i}-a_{j}\right)
\left(  b_{i}-b_{j}\right)  \right)  \\=\prod_{1\leq i<j\leq n}\left(  \left(
a_{j}-a_{i}\right)  \left(  b_{j}-b_{i}\right)  \right)  \\\text{(here, we
have renamed the index }\left(  j,i\right)  \text{ as }\left(  i,j\right)
\\\text{in the product)}}}\\
&  =\mathbf{m}\left(  0,1,\ldots,n-1\right)  \cdot\prod_{1\leq i<j\leq
n}\left(  \underbrace{\left(  a_{j}-a_{i}\right)  \left(  b_{j}-b_{i}\right)
}_{=\left(  a_{i}-a_{j}\right)  \left(  b_{i}-b_{j}\right)  }\right) \\
&  =\mathbf{m}\left(  0,1,\ldots,n-1\right)  \cdot\prod_{1\leq i<j\leq
n}\left(  \left(  a_{i}-a_{j}\right)  \left(  b_{i}-b_{j}\right)  \right)  .
\end{align*}
Thus,%
\begin{align*}
&  \sum_{\sigma\in S_{n}}\left(  -1\right)  ^{\sigma}\left(  \sum_{i=1}%
^{n}a_{i}b_{\sigma\left(  i\right)  }\right)  ^{m}\\
&  =\sum_{\substack{\left(  g_{1},g_{2},\ldots,g_{n}\right)  \in\mathbb{N}%
^{n};\\g_{1}<g_{2}<\cdots<g_{n};\\g_{1}+g_{2}+\cdots+g_{n}=m}}\mathbf{m}%
\left(  g_{1},g_{2},\ldots,g_{n}\right)  \cdot\det\left(  \left(  a_{i}%
^{g_{j}}\right)  _{1\leq i\leq n,\ 1\leq j\leq n}\right)  \cdot\det\left(
\left(  b_{i}^{g_{j}}\right)  _{1\leq i\leq n,\ 1\leq j\leq n}\right) \\
&  =\mathbf{m}\left(  0,1,\ldots,n-1\right)  \cdot\prod_{1\leq i<j\leq
n}\left(  \left(  a_{i}-a_{j}\right)  \left(  b_{i}-b_{j}\right)  \right)  .
\end{align*}
This proves Lemma \ref{lem.sol.powerdet.r1.b}.
\end{proof}

\begin{proof}
[Solution to Exercise \ref{exe.powerdet.r1}.]Exercise \ref{exe.powerdet.r1}
\textbf{(a)} follows from Lemma \ref{lem.sol.powerdet.r1.a}. Exercise
\ref{exe.powerdet.r1} \textbf{(b)} follows from Lemma
\ref{lem.sol.powerdet.r1.b}. Exercise \ref{exe.powerdet.r1} \textbf{(c)}
follows from Lemma \ref{lem.sol.powerdet.r1.c}.
\end{proof}

\subsubsection{Additional observations}

Let us also record two corollaries of Lemma \ref{lem.sol.powerdet.r1.0}:

\begin{corollary}
\label{cor.sol.powerdet.r1.0dist}Let $n\in\mathbb{N}$. Let $\left(
p_{1},p_{2},\ldots,p_{n}\right)  \in\mathbb{N}^{n}$ be such that $p_{1}%
,p_{2},\ldots,p_{n}$ are distinct. Then, $p_{1}+p_{2}+\cdots+p_{n}\geq
\dbinom{n}{2}$.
\end{corollary}

\begin{proof}
[Proof of Corollary \ref{cor.sol.powerdet.r1.0dist}.]Define the subsets
$\mathbf{E}$ and $\mathbf{I}$ of $\mathbb{N}^{n}$ as in Lemma
\ref{lem.sol.powerdet.r1.EI}. Then, Lemma \ref{lem.sol.powerdet.r1.EI} shows
that the map%
\begin{align*}
\mathbf{I}\times S_{n}  &  \rightarrow\mathbf{E},\\
\left(  \left(  g_{1},g_{2},\ldots,g_{n}\right)  ,\sigma\right)   &
\mapsto\left(  g_{\sigma\left(  1\right)  },g_{\sigma\left(  2\right)
},\ldots,g_{\sigma\left(  n\right)  }\right)
\end{align*}
is well-defined and is a bijection. Denote this map by $\Phi$.

Notice that%
\begin{equation}
\mathbf{E}=\left\{  \left(  k_{1},k_{2},\ldots,k_{n}\right)  \in\mathbb{N}%
^{n}\ \mid\ \text{the integers }k_{1},k_{2},\ldots,k_{n}\text{ are
distinct}\right\}  \label{pf.cor.sol.powerdet.r1.0dist.E=}%
\end{equation}
(by the definition of $\mathbf{E}$).

\begin{vershort}
Now, $\left(  p_{1},p_{2},\ldots,p_{n}\right)  $ is an element of
$\mathbb{N}^{n}$ having the property that the integers $p_{1},p_{2}%
,\ldots,p_{n}$ are distinct. In other words,%
\[
\left(  p_{1},p_{2},\ldots,p_{n}\right)  \in\left\{  \left(  k_{1}%
,k_{2},\ldots,k_{n}\right)  \in\mathbb{N}^{n}\ \mid\ \text{the integers }%
k_{1},k_{2},\ldots,k_{n}\text{ are distinct}\right\}  .
\]
In view of (\ref{pf.cor.sol.powerdet.r1.0dist.E=}), this rewrites as $\left(
p_{1},p_{2},\ldots,p_{n}\right)  \in\mathbf{E}$.
\end{vershort}

\begin{verlong}
Now, $\left(  p_{1},p_{2},\ldots,p_{n}\right)  $ is an element of
$\mathbb{N}^{n}$ having the property that the integers $p_{1},p_{2}%
,\ldots,p_{n}$ are distinct. In other words, $\left(  p_{1},p_{2},\ldots
,p_{n}\right)  $ is a $\left(  k_{1},k_{2},\ldots,k_{n}\right)  \in
\mathbb{N}^{n}$ such that the integers $k_{1},k_{2},\ldots,k_{n}$ are
distinct. In other words,%
\[
\left(  p_{1},p_{2},\ldots,p_{n}\right)  \in\left\{  \left(  k_{1}%
,k_{2},\ldots,k_{n}\right)  \in\mathbb{N}^{n}\ \mid\ \text{the integers }%
k_{1},k_{2},\ldots,k_{n}\text{ are distinct}\right\}  .
\]
In view of (\ref{pf.cor.sol.powerdet.r1.0dist.E=}), this rewrites as $\left(
p_{1},p_{2},\ldots,p_{n}\right)  \in\mathbf{E}$.
\end{verlong}

\begin{vershort}
But the map $\Phi$ is surjective (since $\Phi$ is a bijection). In other
words, $\mathbf{E}=\Phi\left(  \mathbf{I}\times S_{n}\right)  $. Thus,
$\left(  p_{1},p_{2},\ldots,p_{n}\right)  \in\mathbf{E}=\Phi\left(
\mathbf{I}\times S_{n}\right)  $. In other words, there exists some
\newline$\left(  \left(  k_{1},k_{2},\ldots,k_{n}\right)  ,\sigma\right)
\in\mathbf{I}\times S_{n}$ such that $\left(  p_{1},p_{2},\ldots,p_{n}\right)
=\Phi\left(  \left(  \left(  k_{1},k_{2},\ldots,k_{n}\right)  ,\sigma\right)
\right)  $. Consider this $\left(  \left(  k_{1},k_{2},\ldots,k_{n}\right)
,\sigma\right)  $. We have $\left(  k_{1},k_{2},\ldots,k_{n}\right)
\in\mathbf{I}$, and thus $k_{1}<k_{2}<\cdots<k_{n}$ (by the definition of
$\mathbf{I}$).
\end{vershort}

\begin{verlong}
But the map $\Phi$ is bijective (since $\Phi$ is a bijection). Hence, $\Phi$
is surjective. In other words, $\mathbf{E}=\Phi\left(  \mathbf{I}\times
S_{n}\right)  $. Thus, $\left(  p_{1},p_{2},\ldots,p_{n}\right)  \in
\mathbf{E}=\Phi\left(  \mathbf{I}\times S_{n}\right)  $. In other words, there
exists some $\left(  \beta,\sigma\right)  \in\mathbf{I}\times S_{n}$ such that
$\left(  p_{1},p_{2},\ldots,p_{n}\right)  =\Phi\left(  \left(  \beta
,\sigma\right)  \right)  $. Consider this $\left(  \beta,\sigma\right)  $.

We have $\left(  \beta,\sigma\right)  \in\mathbf{I}\times S_{n}$. Thus,
$\beta\in\mathbf{I}$ and $\sigma\in S_{n}$.

We have $\sigma\in S_{n}$. But $S_{n}$ is the set of all permutations of the
set $\left\{  1,2,\ldots,n\right\}  $. In other words, $S_{n}$ is the set of
all bijections $\left\{  1,2,\ldots,n\right\}  \rightarrow\left\{
1,2,\ldots,n\right\}  $ (since the permutations of the set $\left\{
1,2,\ldots,n\right\}  $ are exactly the bijections $\left\{  1,2,\ldots
,n\right\}  \rightarrow\left\{  1,2,\ldots,n\right\}  $). Hence, $\sigma$ is a
bijection $\left\{  1,2,\ldots,n\right\}  \rightarrow\left\{  1,2,\ldots
,n\right\}  $ (since $\sigma\in S_{n}$).

We have%
\[
\beta\in\mathbf{I}=\left\{  \left(  k_{1},k_{2},\ldots,k_{n}\right)
\in\mathbb{N}^{n}\ \mid\ k_{1}<k_{2}<\cdots<k_{n}\right\}
\]
(by the definition of $\mathbf{I}$). Hence, $\beta$ can be written in the form
$\beta=\left(  k_{1},k_{2},\ldots,k_{n}\right)  $ for some $\left(
k_{1},k_{2},\ldots,k_{n}\right)  \in\mathbb{N}^{n}$ satisfying $k_{1}%
<k_{2}<\cdots<k_{n}$. Consider this $\left(  k_{1},k_{2},\ldots,k_{n}\right)
$.
\end{verlong}

\begin{vershort}
We have%
\[
\left(  p_{1},p_{2},\ldots,p_{n}\right)  =\Phi\left(  \left(  \left(
k_{1},k_{2},\ldots,k_{n}\right)  ,\sigma\right)  \right)  =\left(
k_{\sigma\left(  1\right)  },k_{\sigma\left(  2\right)  },\ldots
,k_{\sigma\left(  n\right)  }\right)
\]
(by the definition of the map $\Phi$). Thus,%
\begin{align*}
p_{1}+p_{2}+\cdots+p_{n}  &  =k_{\sigma\left(  1\right)  }+k_{\sigma\left(
2\right)  }+\cdots+k_{\sigma\left(  n\right)  }=k_{1}+k_{2}+\cdots+k_{n}\\
&  \ \ \ \ \ \ \ \ \ \ \left(  \text{by Lemma \ref{lem.sol.powerdet.r1.perm-m}
\textbf{(a)}, applied to }g_{i}=k_{i}\right) \\
&  \geq\dbinom{n}{2}\ \ \ \ \ \ \ \ \ \ \left(  \text{by Lemma
\ref{lem.sol.powerdet.r1.0} \textbf{(a)}}\right)  .
\end{align*}
This proves Corollary \ref{cor.sol.powerdet.r1.0dist}. \qedhere

\end{vershort}

\begin{verlong}
We have%
\begin{align*}
\left(  p_{1},p_{2},\ldots,p_{n}\right)   &  =\Phi\left(  \left(
\underbrace{\beta}_{=\left(  k_{1},k_{2},\ldots,k_{n}\right)  },\sigma\right)
\right)  =\Phi\left(  \left(  \left(  k_{1},k_{2},\ldots,k_{n}\right)
,\sigma\right)  \right) \\
&  =\left(  k_{\sigma\left(  1\right)  },k_{\sigma\left(  2\right)  }%
,\ldots,k_{\sigma\left(  n\right)  }\right)  \ \ \ \ \ \ \ \ \ \ \left(
\text{by the definition of the map }\Phi\right)  .
\end{align*}
Thus,%
\begin{align*}
p_{1}+p_{2}+\cdots+p_{n}  &  =k_{\sigma\left(  1\right)  }+k_{\sigma\left(
2\right)  }+\cdots+k_{\sigma\left(  n\right)  }=\underbrace{\sum_{i=1}^{n}%
}_{=\sum_{i\in\left\{  1,2,\ldots,n\right\}  }}k_{\sigma\left(  i\right)  }\\
&  =\sum_{i\in\left\{  1,2,\ldots,n\right\}  }k_{\sigma\left(  i\right)
}=\underbrace{\sum_{i\in\left\{  1,2,\ldots,n\right\}  }}_{=\sum_{i=1}^{n}%
}k_{i}\\
&  \ \ \ \ \ \ \ \ \ \ \left(
\begin{array}
[c]{c}%
\text{here, we have substituted }i\text{ for }\sigma\left(  i\right)  \text{
in the sum, since}\\
\text{the map }\sigma:\left\{  1,2,\ldots,n\right\}  \rightarrow\left\{
1,2,\ldots,n\right\}  \text{ is a bijection}%
\end{array}
\right) \\
&  =\sum_{i=1}^{n}k_{i}=k_{1}+k_{2}+\cdots+k_{n}\geq\dbinom{n}{2}%
\end{align*}
(by Lemma \ref{lem.sol.powerdet.r1.0} \textbf{(a)}). This proves Corollary
\ref{cor.sol.powerdet.r1.0dist}.
\end{verlong}
\end{proof}

\begin{corollary}
\label{cor.sol.powerdet.r1.1c}Let $n\in\mathbb{N}$. Let $a_{1},a_{2}%
,\ldots,a_{n}$ be $n$ elements of $\mathbb{K}$. Let $b_{1},b_{2},\ldots,b_{n}$
be $n$ elements of $\mathbb{K}$.

Let $\left(  k_{1},k_{2},\ldots,k_{n}\right)  \in\mathbb{N}^{n}$ be such that
$k_{1}+k_{2}+\cdots+k_{n}<\dbinom{n}{2}$. Then,%
\[
\sum_{\sigma\in S_{n}}\left(  -1\right)  ^{\sigma}\prod_{i=1}^{n}\left(
a_{i}b_{\sigma\left(  i\right)  }\right)  ^{k_{i}}=0.
\]

\end{corollary}

\begin{proof}
[Proof of Corollary \ref{cor.sol.powerdet.r1.1c}.]The integers $k_{1}%
,k_{2},\ldots,k_{n}$ are not distinct\footnote{\textit{Proof.} Assume the
contrary. Thus, the integers $k_{1},k_{2},\ldots,k_{n}$ are distinct. Hence,
Corollary \ref{cor.sol.powerdet.r1.0dist} (applied to $k_{i}$ instead of
$p_{i}$) shows that $k_{1}+k_{2}+\cdots+k_{n}\geq\dbinom{n}{2}$. This
contradicts $k_{1}+k_{2}+\cdots+k_{n}<\dbinom{n}{2}$. This contradiction shows
that our assumption was wrong. Qed.}. Hence, Lemma \ref{lem.sol.powerdet.r1.1}
\textbf{(b)} yields $\sum_{\sigma\in S_{n}}\left(  -1\right)  ^{\sigma}%
\prod_{i=1}^{n}\left(  a_{i}b_{\sigma\left(  i\right)  }\right)  ^{k_{i}}=0$.
This proves Corollary \ref{cor.sol.powerdet.r1.1c}.
\end{proof}

\subsection{\label{sect.sol.jacobi-complement}Solution to Exercise
\ref{addexe.jacobi-complement}}

\subsubsection{First solution}

Our first solution to Exercise \ref{addexe.jacobi-complement} (inspired by
\cite[Lemma 9.2.10]{BruRys91}) shall follow the hint given. We are going to
prepare for it by stating several simple lemmas.

First, let us agree on some notations. We are going to use the notations
introduced in Definition \ref{def.submatrix}, in Definition
\ref{def.sect.laplace.notations} and in Definition
\ref{def.sol.addexe.jacobi-complement.Agd} throughout Section
\ref{sect.sol.jacobi-complement}.

Now, we collect some useful lemmas.

\begin{lemma}
\label{lem.sol.addexe.jacobi-complement.Agd.sub}Let $n\in\mathbb{N}$ and
$m\in\mathbb{N}$. Let $A\in\mathbb{K}^{n\times m}$. Let $\gamma\in S_{n}$ and
$\delta\in S_{m}$. Let $i_{1},i_{2},\ldots,i_{u}$ be some elements of
$\left\{  1,2,\ldots,n\right\}  $; let $j_{1},j_{2},\ldots,j_{v}$ be some
elements of $\left\{  1,2,\ldots,m\right\}  $. Then,%
\[
\operatorname*{sub}\nolimits_{i_{1},i_{2},\ldots,i_{u}}^{j_{1},j_{2}%
,\ldots,j_{v}}\left(  A_{\left[  \gamma,\delta\right]  }\right)
=\operatorname*{sub}\nolimits_{\gamma\left(  i_{1}\right)  ,\gamma\left(
i_{2}\right)  ,\ldots,\gamma\left(  i_{u}\right)  }^{\delta\left(
j_{1}\right)  ,\delta\left(  j_{2}\right)  ,\ldots,\delta\left(  j_{v}\right)
}A.
\]

\end{lemma}

\begin{proof}
[Proof of Lemma \ref{lem.sol.addexe.jacobi-complement.Agd.sub}.]Write the
$n\times m$-matrix $A$ in the form $A=\left(  a_{i,j}\right)  _{1\leq i\leq
n,\ 1\leq j\leq m}$. Then,
\begin{equation}
\operatorname*{sub}\nolimits_{\gamma\left(  i_{1}\right)  ,\gamma\left(
i_{2}\right)  ,\ldots,\gamma\left(  i_{u}\right)  }^{\delta\left(
j_{1}\right)  ,\delta\left(  j_{2}\right)  ,\ldots,\delta\left(  j_{v}\right)
}A=\left(  a_{\gamma\left(  i_{x}\right)  ,\delta\left(  j_{y}\right)
}\right)  _{1\leq x\leq u,\ 1\leq y\leq v}
\label{pf.lem.sol.addexe.jacobi-complement.Agd.sub.1}%
\end{equation}
(by the definition of $\operatorname*{sub}\nolimits_{\gamma\left(
i_{1}\right)  ,\gamma\left(  i_{2}\right)  ,\ldots,\gamma\left(  i_{u}\right)
}^{\delta\left(  j_{1}\right)  ,\delta\left(  j_{2}\right)  ,\ldots
,\delta\left(  j_{v}\right)  }A$, since $A=\left(  a_{i,j}\right)  _{1\leq
i\leq n,\ 1\leq j\leq m}$). On the other hand, the definition of $A_{\left[
\gamma,\delta\right]  }$ yields $A_{\left[  \gamma,\delta\right]  }=\left(
a_{\gamma\left(  i\right)  ,\delta\left(  j\right)  }\right)  _{1\leq i\leq
n,\ 1\leq j\leq m}$. Hence, the definition of $\operatorname*{sub}%
\nolimits_{i_{1},i_{2},\ldots,i_{u}}^{j_{1},j_{2},\ldots,j_{v}}\left(
A_{\left[  \gamma,\delta\right]  }\right)  $ yields%
\[
\operatorname*{sub}\nolimits_{i_{1},i_{2},\ldots,i_{u}}^{j_{1},j_{2}%
,\ldots,j_{v}}\left(  A_{\left[  \gamma,\delta\right]  }\right)  =\left(
a_{\gamma\left(  i_{x}\right)  ,\delta\left(  j_{y}\right)  }\right)  _{1\leq
x\leq u,\ 1\leq y\leq v}.
\]
Comparing this with (\ref{pf.lem.sol.addexe.jacobi-complement.Agd.sub.1}), we
obtain $\operatorname*{sub}\nolimits_{i_{1},i_{2},\ldots,i_{u}}^{j_{1}%
,j_{2},\ldots,j_{v}}\left(  A_{\left[  \gamma,\delta\right]  }\right)
=\operatorname*{sub}\nolimits_{\gamma\left(  i_{1}\right)  ,\gamma\left(
i_{2}\right)  ,\ldots,\gamma\left(  i_{u}\right)  }^{\delta\left(
j_{1}\right)  ,\delta\left(  j_{2}\right)  ,\ldots,\delta\left(  j_{v}\right)
}A$. This proves Lemma \ref{lem.sol.addexe.jacobi-complement.Agd.sub}.
\end{proof}

\begin{lemma}
\label{lem.sol.addexe.jacobi-complement.Agd.AB}Let $n\in\mathbb{N}$,
$m\in\mathbb{N}$ and $p\in\mathbb{N}$. Let $\gamma\in S_{n}$, $\delta\in
S_{m}$ and $\varepsilon\in S_{p}$. Let $A\in\mathbb{K}^{n\times m}$ and
$B\in\mathbb{K}^{m\times p}$. Then,%
\[
\left(  AB\right)  _{\left[  \gamma,\varepsilon\right]  }=A_{\left[
\gamma,\delta\right]  }B_{\left[  \delta,\varepsilon\right]  }.
\]

\end{lemma}

\begin{proof}
[Proof of Lemma \ref{lem.sol.addexe.jacobi-complement.Agd.AB}.]Write the
$n\times m$-matrix $A$ in the form $A=\left(  a_{i,j}\right)  _{1\leq i\leq
n,\ 1\leq j\leq m}$. Then, $A_{\left[  \gamma,\delta\right]  }=\left(
a_{\gamma\left(  i\right)  ,\delta\left(  j\right)  }\right)  _{1\leq i\leq
n,\ 1\leq j\leq m}$ (by the definition of $A_{\left[  \gamma,\delta\right]  }$).

Write the $m\times p$-matrix $B$ in the form $B=\left(  b_{i,j}\right)
_{1\leq i\leq m,\ 1\leq j\leq p}$. Then, $B_{\left[  \delta,\varepsilon
\right]  }=\left(  b_{\delta\left(  i\right)  ,\varepsilon\left(  j\right)
}\right)  _{1\leq i\leq m,\ 1\leq j\leq p}$ (by the definition of $B_{\left[
\delta,\varepsilon\right]  }$).

\begin{verlong}
We have $\delta\in S_{m}$. In other words, $\delta$ is a permutation of
$\left\{  1,2,\ldots,m\right\}  $ (since $S_{m}$ is the set of all
permutations of $\left\{  1,2,\ldots,m\right\}  $ (by the definition of
$S_{m}$)). In other words, $\delta$ is a bijective map $\left\{
1,2,\ldots,m\right\}  \rightarrow\left\{  1,2,\ldots,m\right\}  $. In other
words, $\delta:\left\{  1,2,\ldots,m\right\}  \rightarrow\left\{
1,2,\ldots,m\right\}  $ is a bijection.
\end{verlong}

\begin{vershort}
Every $\left(  i,j\right)  \in\left\{  1,2,\ldots,n\right\}  \times\left\{
1,2,\ldots,p\right\}  $ satisfies%
\[
\sum_{k=1}^{m}a_{\gamma\left(  i\right)  ,\delta\left(  k\right)  }%
b_{\delta\left(  k\right)  ,\varepsilon\left(  j\right)  }=\sum_{k=1}%
^{m}a_{\gamma\left(  i\right)  ,k}b_{k,\varepsilon\left(  j\right)  }%
\]
(here, we have substituted $k$ for $\delta\left(  k\right)  $ in the sum,
since the map $\delta:\left\{  1,2,\ldots,m\right\}  \rightarrow\left\{
1,2,\ldots,m\right\}  $ is a bijection). In other words, we have%
\begin{equation}
\left(  \sum_{k=1}^{m}a_{\gamma\left(  i\right)  ,\delta\left(  k\right)
}b_{\delta\left(  k\right)  ,\varepsilon\left(  j\right)  }\right)  _{1\leq
i\leq n,\ 1\leq j\leq p}=\left(  \sum_{k=1}^{m}a_{\gamma\left(  i\right)
,k}b_{k,\varepsilon\left(  j\right)  }\right)  _{1\leq i\leq n,\ 1\leq j\leq
p}. \label{pf.lem.sol.addexe.jacobi-complement.Agd.AB.short.1}%
\end{equation}

The definition of the product of two matrices yields%
\begin{align}
A_{\left[  \gamma,\delta\right]  }B_{\left[  \delta,\varepsilon\right]  }  &
=\left(  \sum_{k=1}^{m}a_{\gamma\left(  i\right)  ,\delta\left(  k\right)
}b_{\delta\left(  k\right)  ,\varepsilon\left(  j\right)  }\right)  _{1\leq
i\leq n,\ 1\leq j\leq p}\nonumber\\
&  \ \ \ \ \ \ \ \ \ \ \left(
\begin{array}
[c]{c}%
\text{since }A_{\left[  \gamma,\delta\right]  }=\left(  a_{\gamma\left(
i\right)  ,\delta\left(  j\right)  }\right)  _{1\leq i\leq n,\ 1\leq j\leq
m}\\
\text{and }B_{\left[  \delta,\varepsilon\right]  }=\left(  b_{\delta\left(
i\right)  ,\varepsilon\left(  j\right)  }\right)  _{1\leq i\leq m,\ 1\leq
j\leq p}%
\end{array}
\right) \nonumber\\
&  =\left(  \sum_{k=1}^{m}a_{\gamma\left(  i\right)  ,k}b_{k,\varepsilon
\left(  j\right)  }\right)  _{1\leq i\leq n,\ 1\leq j\leq p}%
\ \ \ \ \ \ \ \ \ \ \left(  \text{by
(\ref{pf.lem.sol.addexe.jacobi-complement.Agd.AB.short.1})}\right)  .
\label{pf.lem.sol.addexe.jacobi-complement.Agd.AB.short.2}%
\end{align}

On the other hand, the definition of the product of two matrices yields
$AB=\left(  \sum_{k=1}^{m}a_{i,k}b_{k,j}\right)  _{1\leq i\leq n,\ 1\leq j\leq
p}$ (since $A=\left(  a_{i,j}\right)  _{1\leq i\leq n,\ 1\leq j\leq m}$ and
$B=\left(  b_{i,j}\right)  _{1\leq i\leq m,\ 1\leq j\leq p}$). Hence, the
definition of $\left(  AB\right)  _{\left[  \gamma,\varepsilon\right]  }$
yields%
\[
\left(  AB\right)  _{\left[  \gamma,\varepsilon\right]  }=\left(  \sum
_{k=1}^{m}a_{\gamma\left(  i\right)  ,k}b_{k,\varepsilon\left(  j\right)
}\right)  _{1\leq i\leq n,\ 1\leq j\leq p}.
\]
Comparing this with (\ref{pf.lem.sol.addexe.jacobi-complement.Agd.AB.short.2}%
), we obtain $\left(  AB\right)  _{\left[  \gamma,\varepsilon\right]
}=A_{\left[  \gamma,\delta\right]  }B_{\left[  \delta,\varepsilon\right]  }$.
Thus, Lemma \ref{lem.sol.addexe.jacobi-complement.Agd.AB} is proven.
\qedhere

\end{vershort}

\begin{verlong}
Every $\left(  i,j\right)  \in\left\{  1,2,\ldots,n\right\}  \times\left\{
1,2,\ldots,p\right\}  $ satisfies%
\begin{align*}
&  \underbrace{\sum_{k=1}^{m}}_{=\sum_{k\in\left\{  1,2,\ldots,m\right\}  }%
}a_{\gamma\left(  i\right)  ,\delta\left(  k\right)  }b_{\delta\left(
k\right)  ,\varepsilon\left(  j\right)  }\\
&  =\sum_{k\in\left\{  1,2,\ldots,m\right\}  }a_{\gamma\left(  i\right)
,\delta\left(  k\right)  }b_{\delta\left(  k\right)  ,\varepsilon\left(
j\right)  }=\underbrace{\sum_{k\in\left\{  1,2,\ldots,m\right\}  }}%
_{=\sum_{k=1}^{m}}a_{\gamma\left(  i\right)  ,k}b_{k,\varepsilon\left(
j\right)  }\\
&  \ \ \ \ \ \ \ \ \ \ \left(
\begin{array}
[c]{c}%
\text{here, we have substituted }k\text{ for }\delta\left(  k\right)  \text{
in the sum, since}\\
\text{the map }\delta:\left\{  1,2,\ldots,m\right\}  \rightarrow\left\{
1,2,\ldots,m\right\}  \text{ is a bijection}%
\end{array}
\right) \\
&  =\sum_{k=1}^{m}a_{\gamma\left(  i\right)  ,k}b_{k,\varepsilon\left(
j\right)  }.
\end{align*}
In other words, we have%
\begin{equation}
\left(  \sum_{k=1}^{m}a_{\gamma\left(  i\right)  ,\delta\left(  k\right)
}b_{\delta\left(  k\right)  ,\varepsilon\left(  j\right)  }\right)  _{1\leq
i\leq n,\ 1\leq j\leq p}=\left(  \sum_{k=1}^{m}a_{\gamma\left(  i\right)
,k}b_{k,\varepsilon\left(  j\right)  }\right)  _{1\leq i\leq n,\ 1\leq j\leq
p}. \label{pf.lem.sol.addexe.jacobi-complement.Agd.AB.1}%
\end{equation}

The definition of the product of two matrices yields%
\begin{align}
A_{\left[  \gamma,\delta\right]  }B_{\left[  \delta,\varepsilon\right]  }  &
=\left(  \sum_{k=1}^{m}a_{\gamma\left(  i\right)  ,\delta\left(  k\right)
}b_{\delta\left(  k\right)  ,\varepsilon\left(  j\right)  }\right)  _{1\leq
i\leq n,\ 1\leq j\leq p}\nonumber\\
&  \ \ \ \ \ \ \ \ \ \ \left(
\begin{array}
[c]{c}%
\text{since }A_{\left[  \gamma,\delta\right]  }=\left(  a_{\gamma\left(
i\right)  ,\delta\left(  j\right)  }\right)  _{1\leq i\leq n,\ 1\leq j\leq
m}\\
\text{and }B_{\left[  \delta,\varepsilon\right]  }=\left(  b_{\delta\left(
i\right)  ,\varepsilon\left(  j\right)  }\right)  _{1\leq i\leq m,\ 1\leq
j\leq p}%
\end{array}
\right) \nonumber\\
&  =\left(  \sum_{k=1}^{m}a_{\gamma\left(  i\right)  ,k}b_{k,\varepsilon
\left(  j\right)  }\right)  _{1\leq i\leq n,\ 1\leq j\leq p}%
\ \ \ \ \ \ \ \ \ \ \left(  \text{by
(\ref{pf.lem.sol.addexe.jacobi-complement.Agd.AB.1})}\right)  .
\label{pf.lem.sol.addexe.jacobi-complement.Agd.AB.2}%
\end{align}

On the other hand, the definition of the product of two matrices yields
$AB=\left(  \sum_{k=1}^{m}a_{i,k}b_{k,j}\right)  _{1\leq i\leq n,\ 1\leq j\leq
p}$ (since $A=\left(  a_{i,j}\right)  _{1\leq i\leq n,\ 1\leq j\leq m}$ and
$B=\left(  b_{i,j}\right)  _{1\leq i\leq m,\ 1\leq j\leq p}$). Hence, the
definition of $\left(  AB\right)  _{\left[  \gamma,\varepsilon\right]  }$
yields%
\[
\left(  AB\right)  _{\left[  \gamma,\varepsilon\right]  }=\left(  \sum
_{k=1}^{m}a_{\gamma\left(  i\right)  ,k}b_{k,\varepsilon\left(  j\right)
}\right)  _{1\leq i\leq n,\ 1\leq j\leq p}.
\]
Comparing this with (\ref{pf.lem.sol.addexe.jacobi-complement.Agd.AB.2}), we
obtain $\left(  AB\right)  _{\left[  \gamma,\varepsilon\right]  }=A_{\left[
\gamma,\delta\right]  }B_{\left[  \delta,\varepsilon\right]  }$. Thus, Lemma
\ref{lem.sol.addexe.jacobi-complement.Agd.AB} is proven.
\end{verlong}
\end{proof}

\begin{lemma}
\label{lem.sol.addexe.jacobi-complement.Agd.I}Let $n\in\mathbb{N}$ and
$\gamma\in S_{n}$. Then, $\left(  I_{n}\right)  _{\left[  \gamma
,\gamma\right]  }=I_{n}$.
\end{lemma}

\begin{vershort}
\begin{proof}
[Proof of Lemma \ref{lem.sol.addexe.jacobi-complement.Agd.I}.]For every two
objects $i$ and $j$, define an element $\delta_{i,j}\in\mathbb{K}$ by
$\delta_{i,j}=%
\begin{cases}
1, & \text{if }i=j;\\
0, & \text{if }i\neq j
\end{cases}
$. Then, $I_{n}=\left(  \delta_{i,j}\right)  _{1\leq i\leq n,\ 1\leq j\leq n}$
(by the definition of $I_{n}$). Hence, $\left(  I_{n}\right)  _{\left[
\gamma,\gamma\right]  }=\left(  \delta_{\gamma\left(  i\right)  ,\gamma\left(
j\right)  }\right)  _{1\leq i\leq n,\ 1\leq j\leq n}$ (by the definition of
$\left(  I_{n}\right)  _{\left[  \gamma,\gamma\right]  }$).

The map $\gamma$ is a permutation (since $\gamma\in S_{n}$), thus bijective,
thus injective.

For every $\left(  i,j\right)  \in\left\{  1,2,\ldots,n\right\}  ^{2}$, we
have $\delta_{\gamma\left(  i\right)  ,\gamma\left(  j\right)  }=\delta_{i,j}%
$\ \ \ \ \footnote{\textit{Proof.} Let $\left(  i,j\right)  \in\left\{
1,2,\ldots,n\right\}  ^{2}$. We must prove $\delta_{\gamma\left(  i\right)
,\gamma\left(  j\right)  }=\delta_{i,j}$.
\par
If $i=j$, then both $\delta_{\gamma\left(  i\right)  ,\gamma\left(  j\right)
}$ and $\delta_{i,j}$ equal $1$ (because $\gamma\left(  \underbrace{i}%
_{=j}\right)  =\gamma\left(  j\right)  $ shows that $\delta_{\gamma\left(
i\right)  ,\gamma\left(  j\right)  }=1$, whereas $i=j$ shows that
$\delta_{i,j}=1$). Hence, if $i=j$, then $\delta_{\gamma\left(  i\right)
,\gamma\left(  j\right)  }=\delta_{i,j}$ holds. Thus, for the rest of our
proof of $\delta_{\gamma\left(  i\right)  ,\gamma\left(  j\right)  }%
=\delta_{i,j}$, we WLOG assume that $i\neq j$. Hence, $\delta_{i,j}=0$.
\par
But the map $\gamma$ is injective. Thus, from $i\neq j$, we obtain
$\gamma\left(  i\right)  \neq\gamma\left(  j\right)  $. Hence, $\delta
_{\gamma\left(  i\right)  ,\gamma\left(  j\right)  }=0$. Comparing this with
$\delta_{i,j}=0$, we obtain $\delta_{\gamma\left(  i\right)  ,\gamma\left(
j\right)  }=\delta_{i,j}$, qed.}. In other words, $\left(  \delta
_{\gamma\left(  i\right)  ,\gamma\left(  j\right)  }\right)  _{1\leq i\leq
n,\ 1\leq j\leq n}=\left(  \delta_{i,j}\right)  _{1\leq i\leq n,\ 1\leq j\leq
n}$. Hence,
\[
\left(  I_{n}\right)  _{\left[  \gamma,\gamma\right]  }=\left(  \delta
_{\gamma\left(  i\right)  ,\gamma\left(  j\right)  }\right)  _{1\leq i\leq
n,\ 1\leq j\leq n}=\left(  \delta_{i,j}\right)  _{1\leq i\leq n,\ 1\leq j\leq
n}.
\]
Comparing this with $I_{n}=\left(  \delta_{i,j}\right)  _{1\leq i\leq
n,\ 1\leq j\leq n}$, we obtain $\left(  I_{n}\right)  _{\left[  \gamma
,\gamma\right]  }=I_{n}$. This proves Lemma
\ref{lem.sol.addexe.jacobi-complement.Agd.I}.
\end{proof}
\end{vershort}

\begin{verlong}
\begin{proof}
[Proof of Lemma \ref{lem.sol.addexe.jacobi-complement.Agd.I}.]For every two
objects $i$ and $j$, define an element $\delta_{i,j}\in\mathbb{K}$ by
$\delta_{i,j}=%
\begin{cases}
1, & \text{if }i=j;\\
0, & \text{if }i\neq j
\end{cases}
$. Then, $I_{n}=\left(  \delta_{i,j}\right)  _{1\leq i\leq n,\ 1\leq j\leq n}$
(by the definition of $I_{n}$). Hence, $\left(  I_{n}\right)  _{\left[
\gamma,\gamma\right]  }=\left(  \delta_{\gamma\left(  i\right)  ,\gamma\left(
j\right)  }\right)  _{1\leq i\leq n,\ 1\leq j\leq n}$ (by the definition of
$\left(  I_{n}\right)  _{\left[  \gamma,\gamma\right]  }$).

We have $\gamma\in S_{n}$. In other words, $\gamma$ is a permutation of
$\left\{  1,2,\ldots,n\right\}  $ (since $S_{n}$ is the set of all
permutations of $\left\{  1,2,\ldots,n\right\}  $ (by the definition of
$S_{n}$)). In other words, $\gamma$ is a bijective map $\left\{
1,2,\ldots,n\right\}  \rightarrow\left\{  1,2,\ldots,n\right\}  $. The map
$\gamma$ is bijective, and thus injective.

For every $\left(  i,j\right)  \in\left\{  1,2,\ldots,n\right\}  ^{2}$, we
have%
\begin{equation}
\delta_{\gamma\left(  i\right)  ,\gamma\left(  j\right)  }=\delta_{i,j}
\label{pf.lem.sol.addexe.jacobi-complement.Agd.I.1}%
\end{equation}
\footnote{\textit{Proof of (\ref{pf.lem.sol.addexe.jacobi-complement.Agd.I.1}%
):} Let $\left(  i,j\right)  \in\left\{  1,2,\ldots,n\right\}  ^{2}$. We must
prove (\ref{pf.lem.sol.addexe.jacobi-complement.Agd.I.1}).
\par
We have $\left(  i,j\right)  \in\left\{  1,2,\ldots,n\right\}  ^{2}$. In other
words, $i\in\left\{  1,2,\ldots,n\right\}  $ and $j\in\left\{  1,2,\ldots
,n\right\}  $. We are in one of the following two cases:
\par
\textit{Case 1:} We have $i=j$.
\par
\textit{Case 2:} We have $i\neq j$.
\par
Let us first consider Case 1. In this case, we have $i=j$. Thus, $\delta
_{i,j}=1$. But $\gamma\left(  \underbrace{i}_{=j}\right)  =\gamma\left(
j\right)  $, and thus $\delta_{\gamma\left(  i\right)  ,\gamma\left(
j\right)  }=1$. Comparing this with $\delta_{i,j}=1$, we obtain $\delta
_{\gamma\left(  i\right)  ,\gamma\left(  j\right)  }=\delta_{i,j}$. Thus,
(\ref{pf.lem.sol.addexe.jacobi-complement.Agd.I.1}) is proven in Case 1.
\par
Let us now consider Case 2. In this case, we have $i\neq j$. In other words,
the elements $i$ and $j$ are distinct. But the map $\gamma$ is injective.
Hence, if $u$ and $v$ are two distinct elements of $\left\{  1,2,\ldots
,n\right\}  $, then $\gamma\left(  u\right)  \neq\gamma\left(  v\right)  $.
Applying this to $u=i$ and $v=j$, we conclude that $\gamma\left(  i\right)
\neq\gamma\left(  j\right)  $ (since $i$ and $j$ are distinct). Hence,
$\delta_{\gamma\left(  i\right)  ,\gamma\left(  j\right)  }=0$. Comparing this
with $\delta_{i,j}=0$ (since $i\neq j$), we obtain $\delta_{\gamma\left(
i\right)  ,\gamma\left(  j\right)  }=\delta_{i,j}$. Thus,
(\ref{pf.lem.sol.addexe.jacobi-complement.Agd.I.1}) is proven in Case 2.
\par
We have now proven (\ref{pf.lem.sol.addexe.jacobi-complement.Agd.I.1}) in each
of the two Cases 1 and 2. Since these two Cases cover all possibilities, we
can thus conclude that (\ref{pf.lem.sol.addexe.jacobi-complement.Agd.I.1})
always holds. Thus, (\ref{pf.lem.sol.addexe.jacobi-complement.Agd.I.1}) is
proven.}. In other words, $\left(  \delta_{\gamma\left(  i\right)
,\gamma\left(  j\right)  }\right)  _{1\leq i\leq n,\ 1\leq j\leq n}=\left(
\delta_{i,j}\right)  _{1\leq i\leq n,\ 1\leq j\leq n}$. Hence,
\[
\left(  I_{n}\right)  _{\left[  \gamma,\gamma\right]  }=\left(  \delta
_{\gamma\left(  i\right)  ,\gamma\left(  j\right)  }\right)  _{1\leq i\leq
n,\ 1\leq j\leq n}=\left(  \delta_{i,j}\right)  _{1\leq i\leq n,\ 1\leq j\leq
n}.
\]
Comparing this with $I_{n}=\left(  \delta_{i,j}\right)  _{1\leq i\leq
n,\ 1\leq j\leq n}$, we obtain $\left(  I_{n}\right)  _{\left[  \gamma
,\gamma\right]  }=I_{n}$. This proves Lemma
\ref{lem.sol.addexe.jacobi-complement.Agd.I}.
\end{proof}
\end{verlong}

\begin{lemma}
\label{lem.sol.addexe.jacobi-complement.Agd.inv}Let $n\in\mathbb{N}$. Let
$A\in\mathbb{K}^{n\times n}$ be an invertible matrix. Let $\gamma\in S_{n}$
and $\delta\in S_{n}$. Then, the matrix $A_{\left[  \gamma,\delta\right]  }%
\in\mathbb{K}^{n\times n}$ is invertible, and its inverse is $\left(
A_{\left[  \gamma,\delta\right]  }\right)  ^{-1}=\left(  A^{-1}\right)
_{\left[  \delta,\gamma\right]  }$.
\end{lemma}

\begin{proof}
[Proof of Lemma \ref{lem.sol.addexe.jacobi-complement.Agd.inv}.]Lemma
\ref{lem.sol.addexe.jacobi-complement.Agd.AB} (applied to $n$, $n$, $\gamma$
and $A^{-1}$ instead of $m$, $p$, $\varepsilon$ and $B$) yields $\left(
AA^{-1}\right)  _{\left[  \gamma,\gamma\right]  }=A_{\left[  \gamma
,\delta\right]  }\left(  A^{-1}\right)  _{\left[  \delta,\gamma\right]  }$.
Thus,%
\begin{equation}
A_{\left[  \gamma,\delta\right]  }\left(  A^{-1}\right)  _{\left[
\delta,\gamma\right]  }=\left(  \underbrace{AA^{-1}}_{=I_{n}}\right)
_{\left[  \gamma,\gamma\right]  }=\left(  I_{n}\right)  _{\left[
\gamma,\gamma\right]  }=I_{n}
\label{pf.lem.sol.addexe.jacobi-complement.Agd.inv.1}%
\end{equation}
(by Lemma \ref{lem.sol.addexe.jacobi-complement.Agd.I}). Also, Lemma
\ref{lem.sol.addexe.jacobi-complement.Agd.AB} (applied to $n$, $n$, $\delta$,
$\gamma$, $\delta$, $A^{-1}$ and $A$ instead of $m$, $p$, $\gamma$, $\delta$,
$\varepsilon$, $A$ and $B$) yields $\left(  A^{-1}A\right)  _{\left[
\delta,\delta\right]  }=\left(  A^{-1}\right)  _{\left[  \delta,\gamma\right]
}A_{\left[  \gamma,\delta\right]  }$. Thus,%
\begin{equation}
\left(  A^{-1}\right)  _{\left[  \delta,\gamma\right]  }A_{\left[
\gamma,\delta\right]  }=\left(  \underbrace{A^{-1}A}_{=I_{n}}\right)
_{\left[  \delta,\delta\right]  }=\left(  I_{n}\right)  _{\left[
\delta,\delta\right]  }=I_{n}
\label{pf.lem.sol.addexe.jacobi-complement.Agd.inv.2}%
\end{equation}
(by Lemma \ref{lem.sol.addexe.jacobi-complement.Agd.I} (applied to $\delta$
instead of $\gamma$)).

\begin{vershort}
The two equalities (\ref{pf.lem.sol.addexe.jacobi-complement.Agd.inv.1}) and
(\ref{pf.lem.sol.addexe.jacobi-complement.Agd.inv.2}) (combined) show that the
matrix $\left(  A^{-1}\right)  _{\left[  \delta,\gamma\right]  }$ is an
inverse of the matrix $A_{\left[  \gamma,\delta\right]  }$. Thus, the matrix
$A_{\left[  \gamma,\delta\right]  }\in\mathbb{K}^{n\times n}$ is invertible,
and its inverse is $\left(  A_{\left[  \gamma,\delta\right]  }\right)
^{-1}=\left(  A^{-1}\right)  _{\left[  \delta,\gamma\right]  }$. Lemma
\ref{lem.sol.addexe.jacobi-complement.Agd.inv} is proven. \qedhere

\end{vershort}

\begin{verlong}
So we have shown that the $n\times n$-matrix $\left(  A^{-1}\right)  _{\left[
\delta,\gamma\right]  }$ satisfies $\left(  A^{-1}\right)  _{\left[
\delta,\gamma\right]  }A_{\left[  \gamma,\delta\right]  }=I_{n}$ and
$A_{\left[  \gamma,\delta\right]  }\left(  A^{-1}\right)  _{\left[
\delta,\gamma\right]  }=I_{n}$. In other words, $\left(  A^{-1}\right)
_{\left[  \delta,\gamma\right]  }$ is an $n\times n$-matrix $B$ such that
$BA_{\left[  \gamma,\delta\right]  }=I_{n}$ and $A_{\left[  \gamma
,\delta\right]  }B=I_{n}$. In other words, $\left(  A^{-1}\right)  _{\left[
\delta,\gamma\right]  }$ is an inverse of $A_{\left[  \gamma,\delta\right]  }$
(because an inverse of $A_{\left[  \gamma,\delta\right]  }$ means an $n\times
n$-matrix $B$ such that $BA_{\left[  \gamma,\delta\right]  }=I_{n}$ and
$A_{\left[  \gamma,\delta\right]  }B=I_{n}$ (by the definition of an
\textquotedblleft inverse of $A_{\left[  \gamma,\delta\right]  }%
$\textquotedblright)). Thus, the matrix $A_{\left[  \gamma,\delta\right]  }$
has an inverse; in other words, the matrix $A_{\left[  \gamma,\delta\right]
}$ is invertible. Also, clearly, the inverse of $A_{\left[  \gamma
,\delta\right]  }$ is $\left(  A^{-1}\right)  _{\left[  \delta,\gamma\right]
}$ (since we have proven that $\left(  A^{-1}\right)  _{\left[  \delta
,\gamma\right]  }$ is an inverse of $A_{\left[  \gamma,\delta\right]  }$). In
other words, $\left(  A_{\left[  \gamma,\delta\right]  }\right)  ^{-1}=\left(
A^{-1}\right)  _{\left[  \delta,\gamma\right]  }$. Thus, the inverse of
$A_{\left[  \gamma,\delta\right]  }$ is $\left(  A_{\left[  \gamma
,\delta\right]  }\right)  ^{-1}=\left(  A^{-1}\right)  _{\left[  \delta
,\gamma\right]  }$. This completes the proof of Lemma
\ref{lem.sol.addexe.jacobi-complement.Agd.inv}.
\end{verlong}
\end{proof}

We can finally step to the solution of Exercise \ref{addexe.jacobi-complement}:

\begin{proof}
[First solution to Exercise \ref{addexe.jacobi-complement}.]Let $k=\left\vert
P\right\vert $. Thus, $k=\left\vert P\right\vert =\left\vert Q\right\vert $.

\begin{vershort}
Clearly, $k\in\left\{  0,1,\ldots,n\right\}  $ (since $k=\left\vert
P\right\vert $ for a subset $P$ of $\left\{  1,2,\ldots,n\right\}  $).
\end{vershort}

\begin{verlong}
We know that $P$ is a subset of $\left\{  1,2,\ldots,n\right\}  $. Thus, $P$
is a finite set (since $\left\{  1,2,\ldots,n\right\}  $ is a finite set).
Hence, $\left\vert P\right\vert \in\mathbb{N}$, so that $k=\left\vert
P\right\vert \in\mathbb{N}$.

The definition of $\widetilde{P}$ yields $\widetilde{P}=\left\{
1,2,\ldots,n\right\}  \setminus P\subseteq\left\{  1,2,\ldots,n\right\}  $.
Thus, $\widetilde{P}$ is a finite set (since $\left\{  1,2,\ldots,n\right\}  $
is a finite set). Hence, $\left\vert \widetilde{P}\right\vert \in\mathbb{N}$.

Also,%
\begin{align*}
\left\vert \underbrace{\widetilde{P}}_{=\left\{  1,2,\ldots,n\right\}
\setminus P}\right\vert  &  =\left\vert \left\{  1,2,\ldots,n\right\}
\setminus P\right\vert =\underbrace{\left\vert \left\{  1,2,\ldots,n\right\}
\right\vert }_{=n}-\underbrace{\left\vert P\right\vert }_{=k}%
\ \ \ \ \ \ \ \ \ \ \left(  \text{since }P\subseteq\left\{  1,2,\ldots
,n\right\}  \right) \\
&  =n-k.
\end{align*}
Hence, $n-k=\left\vert \widetilde{P}\right\vert \in\mathbb{N}$; thus,
$n-k\geq0$. In other words, $k\leq n$. Combined with $k\geq0$ (since
$k=\left\vert P\right\vert \in\mathbb{N}$), this yields $k\in\left\{
0,1,\ldots,n\right\}  $.
\end{verlong}

Lemma \ref{lem.sol.addexe.jacobi-complement.Ialbe} (applied to $I=P$) yields
that there exists a $\sigma\in S_{n}$ satisfying $\left(  \sigma\left(
1\right)  ,\sigma\left(  2\right)  ,\ldots,\sigma\left(  k\right)  \right)
=w\left(  P\right)  $, $\left(  \sigma\left(  k+1\right)  ,\sigma\left(
k+2\right)  ,\ldots,\sigma\left(  n\right)  \right)  =w\left(  \widetilde{P}%
\right)  $ and $\left(  -1\right)  ^{\sigma}=\left(  -1\right)  ^{\sum
P-\left(  1+2+\cdots+k\right)  }$. Denote this $\sigma$ by $\gamma$. Thus,
$\gamma$ is an element of $S_{n}$ satisfying $\left(  \gamma\left(  1\right)
,\gamma\left(  2\right)  ,\ldots,\gamma\left(  k\right)  \right)  =w\left(
P\right)  $, $\left(  \gamma\left(  k+1\right)  ,\gamma\left(  k+2\right)
,\ldots,\gamma\left(  n\right)  \right)  =w\left(  \widetilde{P}\right)  $ and
$\left(  -1\right)  ^{\gamma}=\left(  -1\right)  ^{\sum P-\left(
1+2+\cdots+k\right)  }$.

Lemma \ref{lem.sol.addexe.jacobi-complement.Ialbe} (applied to $I=Q$) yields
that there exists a $\sigma\in S_{n}$ satisfying $\left(  \sigma\left(
1\right)  ,\sigma\left(  2\right)  ,\ldots,\sigma\left(  k\right)  \right)
=w\left(  Q\right)  $, $\left(  \sigma\left(  k+1\right)  ,\sigma\left(
k+2\right)  ,\ldots,\sigma\left(  n\right)  \right)  =w\left(  \widetilde{Q}%
\right)  $ and $\left(  -1\right)  ^{\sigma}=\left(  -1\right)  ^{\sum
Q-\left(  1+2+\cdots+k\right)  }$. Denote this $\sigma$ by $\delta$. Thus,
$\delta$ is an element of $S_{n}$ satisfying $\left(  \delta\left(  1\right)
,\delta\left(  2\right)  ,\ldots,\delta\left(  k\right)  \right)  =w\left(
Q\right)  $, $\left(  \delta\left(  k+1\right)  ,\delta\left(  k+2\right)
,\ldots,\delta\left(  n\right)  \right)  =w\left(  \widetilde{Q}\right)  $ and
$\left(  -1\right)  ^{\delta}=\left(  -1\right)  ^{\sum Q-\left(
1+2+\cdots+k\right)  }$.

Lemma \ref{lem.sol.addexe.jacobi-complement.Agd.inv} shows that the matrix
$A_{\left[  \gamma,\delta\right]  }\in\mathbb{K}^{n\times n}$ is invertible,
and its inverse is $\left(  A_{\left[  \gamma,\delta\right]  }\right)
^{-1}=\left(  A^{-1}\right)  _{\left[  \delta,\gamma\right]  }$. Hence,
Exercise \ref{exe.block2x2.jacobi.rewr} (applied to $A_{\left[  \gamma
,\delta\right]  }$ instead of $A$) yields that%
\begin{align}
&  \det\left(  \operatorname*{sub}\nolimits_{1,2,\ldots,k}^{1,2,\ldots
,k}\left(  A_{\left[  \gamma,\delta\right]  }\right)  \right) \nonumber\\
&  =\det\left(  A_{\left[  \gamma,\delta\right]  }\right)  \cdot\det\left(
\operatorname*{sub}\nolimits_{k+1,k+2,\ldots,n}^{k+1,k+2,\ldots,n}\left(
\underbrace{\left(  A_{\left[  \gamma,\delta\right]  }\right)  ^{-1}%
}_{=\left(  A^{-1}\right)  _{\left[  \delta,\gamma\right]  }}\right)  \right)
\nonumber\\
&  =\det\left(  A_{\left[  \gamma,\delta\right]  }\right)  \cdot\det\left(
\operatorname*{sub}\nolimits_{k+1,k+2,\ldots,n}^{k+1,k+2,\ldots,n}\left(
\left(  A^{-1}\right)  _{\left[  \delta,\gamma\right]  }\right)  \right)  .
\label{sol.addexe.jacobi-complement.1}%
\end{align}
But Lemma \ref{lem.sol.addexe.jacobi-complement.Agd.sub} (applied to $n$, $k$,
$k$, $\left(  1,2,\ldots,k\right)  $ and $\left(  1,2,\ldots,k\right)  $
instead of $m$, $u$, $v$, $\left(  i_{1},i_{2},\ldots,i_{u}\right)  $ and
$\left(  j_{1},j_{2},\ldots,j_{v}\right)  $) yields%
\begin{align}
\operatorname*{sub}\nolimits_{1,2,\ldots,k}^{1,2,\ldots,k}\left(  A_{\left[
\gamma,\delta\right]  }\right)   &  =\operatorname*{sub}\nolimits_{\gamma
\left(  1\right)  ,\gamma\left(  2\right)  ,\ldots,\gamma\left(  k\right)
}^{\delta\left(  1\right)  ,\delta\left(  2\right)  ,\ldots,\delta\left(
k\right)  }A=\operatorname*{sub}\nolimits_{\left(  \gamma\left(  1\right)
,\gamma\left(  2\right)  ,\ldots,\gamma\left(  k\right)  \right)  }^{\left(
\delta\left(  1\right)  ,\delta\left(  2\right)  ,\ldots,\delta\left(
k\right)  \right)  }A\nonumber\\
&  =\operatorname*{sub}\nolimits_{w\left(  P\right)  }^{w\left(  Q\right)  }A
\label{sol.addexe.jacobi-complement.3a}%
\end{align}
(since $\left(  \gamma\left(  1\right)  ,\gamma\left(  2\right)
,\ldots,\gamma\left(  k\right)  \right)  =w\left(  P\right)  $ and $\left(
\delta\left(  1\right)  ,\delta\left(  2\right)  ,\ldots,\delta\left(
k\right)  \right)  =w\left(  Q\right)  $).

Furthermore, we have $A^{-1}\in\mathbb{K}^{n\times n}$ (since $A\in
\mathbb{K}^{n\times n}$). Hence, Lemma
\ref{lem.sol.addexe.jacobi-complement.Agd.sub} (applied to $n$, $A^{-1}$,
$\delta$, $\gamma$, $n-k$, $n-k$, $\left(  k+1,k+2,\ldots,n\right)  $ and
$\left(  k+1,k+2,\ldots,n\right)  $ instead of $m$, $A$, $\gamma$, $\delta$,
$u$, $v$, $\left(  i_{1},i_{2},\ldots,i_{u}\right)  $ and $\left(  j_{1}%
,j_{2},\ldots,j_{v}\right)  $) yields%
\begin{align}
\operatorname*{sub}\nolimits_{k+1,k+2,\ldots,n}^{k+1,k+2,\ldots,n}\left(
\left(  A^{-1}\right)  _{\left[  \delta,\gamma\right]  }\right)   &
=\operatorname*{sub}\nolimits_{\delta\left(  k+1\right)  ,\delta\left(
k+2\right)  ,\ldots,\delta\left(  n\right)  }^{\gamma\left(  k+1\right)
,\gamma\left(  k+2\right)  ,\ldots,\gamma\left(  n\right)  }\left(
A^{-1}\right) \nonumber\\
&  =\operatorname*{sub}\nolimits_{\left(  \delta\left(  k+1\right)
,\delta\left(  k+2\right)  ,\ldots,\delta\left(  n\right)  \right)  }^{\left(
\gamma\left(  k+1\right)  ,\gamma\left(  k+2\right)  ,\ldots,\gamma\left(
n\right)  \right)  }\left(  A^{-1}\right) \nonumber\\
&  =\operatorname*{sub}\nolimits_{w\left(  \widetilde{Q}\right)  }^{w\left(
\widetilde{P}\right)  }\left(  A^{-1}\right)
\label{sol.addexe.jacobi-complement.3b}%
\end{align}
(since $\left(  \delta\left(  k+1\right)  ,\delta\left(  k+2\right)
,\ldots,\delta\left(  n\right)  \right)  =w\left(  \widetilde{Q}\right)  $ and
$\left(  \gamma\left(  k+1\right)  ,\gamma\left(  k+2\right)  ,\ldots
,\gamma\left(  n\right)  \right)  =w\left(  \widetilde{P}\right)  $).

On the other hand,%
\begin{align}
&  \underbrace{\left(  -1\right)  ^{\gamma}}_{=\left(  -1\right)  ^{\sum
P-\left(  1+2+\cdots+k\right)  }}\underbrace{\left(  -1\right)  ^{\delta}%
}_{=\left(  -1\right)  ^{\sum Q-\left(  1+2+\cdots+k\right)  }}\nonumber\\
&  =\left(  -1\right)  ^{\sum P-\left(  1+2+\cdots+k\right)  }\left(
-1\right)  ^{\sum Q-\left(  1+2+\cdots+k\right)  }=\left(  -1\right)
^{\left(  \sum P-\left(  1+2+\cdots+k\right)  \right)  +\left(  \sum Q-\left(
1+2+\cdots+k\right)  \right)  }\nonumber\\
&  =\left(  -1\right)  ^{\sum P+\sum Q} \label{sol.addexe.jacobi-complement.2}%
\end{align}
(since
\begin{align*}
&  \left(  \sum P-\left(  1+2+\cdots+k\right)  \right)  +\left(  \sum
Q-\left(  1+2+\cdots+k\right)  \right) \\
&  =\left(  \sum P+\sum Q\right)  -2\left(  1+2+\cdots+k\right)  \equiv\sum
P+\sum Q\operatorname{mod}2
\end{align*}
).

Lemma \ref{lem.sol.addexe.jacobi-complement.Agd.det} yields%
\begin{equation}
\det\left(  A_{\left[  \gamma,\delta\right]  }\right)  =\underbrace{\left(
-1\right)  ^{\gamma}\left(  -1\right)  ^{\delta}}_{\substack{=\left(
-1\right)  ^{\sum P+\sum Q}\\\text{(by (\ref{sol.addexe.jacobi-complement.2}%
))}}}\det A=\left(  -1\right)  ^{\sum P+\sum Q}\det A.
\label{sol.addexe.jacobi-complement.3c}%
\end{equation}

Now,%
\begin{align*}
&  \det\left(  \underbrace{\operatorname*{sub}\nolimits_{w\left(  P\right)
}^{w\left(  Q\right)  }A}_{\substack{=\operatorname*{sub}\nolimits_{1,2,\ldots
,k}^{1,2,\ldots,k}\left(  A_{\left[  \gamma,\delta\right]  }\right)
\\\text{(by (\ref{sol.addexe.jacobi-complement.3a}))}}}\right) \\
&  =\det\left(  \operatorname*{sub}\nolimits_{1,2,\ldots,k}^{1,2,\ldots
,k}\left(  A_{\left[  \gamma,\delta\right]  }\right)  \right) \\
&  =\underbrace{\det\left(  A_{\left[  \gamma,\delta\right]  }\right)
}_{\substack{=\left(  -1\right)  ^{\sum P+\sum Q}\det A\\\text{(by
(\ref{sol.addexe.jacobi-complement.3c}))}}}\cdot\det\left(
\underbrace{\operatorname*{sub}\nolimits_{k+1,k+2,\ldots,n}^{k+1,k+2,\ldots
,n}\left(  \left(  A^{-1}\right)  _{\left[  \delta,\gamma\right]  }\right)
}_{\substack{=\operatorname*{sub}\nolimits_{w\left(  \widetilde{Q}\right)
}^{w\left(  \widetilde{P}\right)  }\left(  A^{-1}\right)  \\\text{(by
(\ref{sol.addexe.jacobi-complement.3b}))}}}\right)
\ \ \ \ \ \ \ \ \ \ \left(  \text{by (\ref{sol.addexe.jacobi-complement.1}%
)}\right) \\
&  =\left(  -1\right)  ^{\sum P+\sum Q}\det A\cdot\det\left(
\operatorname*{sub}\nolimits_{w\left(  \widetilde{Q}\right)  }^{w\left(
\widetilde{P}\right)  }\left(  A^{-1}\right)  \right)  .
\end{align*}
This solves Exercise \ref{addexe.jacobi-complement}.
\end{proof}

\subsubsection{Second solution}

Now, let us give a second solution to Exercise \ref{addexe.jacobi-complement},
following an idea from \cite[Chapter SCHUR, proof of (1.9)]{LLPT95}.

Throughout this section, we shall use the notations introduced in Definition
\ref{def.submatrix} and in Definition \ref{def.sect.laplace.notations}.

Let us first prove a simple combinatorial fact:

\begin{lemma}
\label{lem.sol.addexe.jacobi-complement.increase-bij}Let $S$ be a set of
integers. Let $u\in\mathbb{N}$. Let $\mathcal{P}_{u}\left(  S\right)  $ denote
the set of all $u$-element subsets of $S$. Let $\mathbf{I}$ denote the set
\[
\left\{  \left(  g_{1},g_{2},\ldots,g_{u}\right)  \in S^{u}\ \mid\ g_{1}%
<g_{2}<\cdots<g_{u}\right\}  .
\]
Then, the map%
\begin{align*}
\mathcal{P}_{u}\left(  S\right)   &  \rightarrow\mathbf{I},\\
R  &  \mapsto w\left(  R\right)
\end{align*}
is well-defined and a bijection.
\end{lemma}

\begin{vershort}
\begin{proof}
[Proof of Lemma \ref{lem.sol.addexe.jacobi-complement.increase-bij}.]Recall
that $\mathcal{P}_{u}\left(  S\right)  $ is the set of all $u$-element subsets
of $S$, whereas $\mathbf{I}$ is the set of all strictly increasing
lists\footnote{A list $\left(  g_{1},g_{2},\ldots,g_{u}\right)  $ is said to
be \textit{strictly increasing} if it satisfies $g_{1}<g_{2}<\cdots<g_{u}$.}
of $u$ elements of $S$.

For every $R\in\mathcal{P}_{u}\left(  S\right)  $, the list $w\left(
R\right)  $ is the list of all elements of $R$ in increasing order (with no
repetitions). Thus, this list $w\left(  R\right)  $ is a strictly increasing
list with $\left\vert R\right\vert =u$ entries (which all belong to $R$ and
therefore to $S$); hence, this list $w\left(  R\right)  $ belongs to
$\mathbf{I}$. Thus, the map
\begin{align*}
\mathcal{P}_{u}\left(  S\right)   &  \rightarrow\mathbf{I},\\
R  &  \mapsto w\left(  R\right)
\end{align*}
is well-defined. Also, the map%
\begin{align*}
\mathbf{I}  &  \rightarrow\mathcal{P}_{u}\left(  S\right)  ,\\
\left(  g_{1},g_{2},\ldots,g_{u}\right)   &  \mapsto\left\{  g_{1}%
,g_{2},\ldots,g_{u}\right\}
\end{align*}
is well-defined (because if $\left(  g_{1},g_{2},\ldots,g_{u}\right)
\in\mathbf{I}$, then we have $g_{1}<g_{2}<\cdots<g_{u}$, and thus the $u$
elements $g_{1},g_{2},\ldots,g_{u}$ are pairwise distinct; therefore,
$\left\{  g_{1},g_{2},\ldots,g_{u}\right\}  $ is a $u$-element set). These two
maps are mutually inverse\footnote{\textit{Proof.} The former map takes a
$u$-element subset $R$ of $S$ and lists its elements in increasing order; the
latter map takes a strictly increasing list (of $u$ elements of $S$) and
outputs the set of its entries. These two operations clearly revert one
another. Thus, the two maps are mutually inverse.}. Hence, the map%
\begin{align*}
\mathcal{P}_{u}\left(  S\right)   &  \rightarrow\mathbf{I},\\
R  &  \mapsto w\left(  R\right)
\end{align*}
is invertible, i.e., is a bijection.
\end{proof}
\end{vershort}

\begin{verlong}
The proof of this lemma relies on the following yet simpler facts:

\begin{lemma}
\label{lem.sol.addexe.jacobi-complement.strinc}Let $u\in\mathbb{N}$. Let
$\varphi:\left\{  1,2,\ldots,u\right\}  \rightarrow\left\{  1,2,\ldots
,u\right\}  $ be a map such that $\varphi\left(  1\right)  <\varphi\left(
2\right)  <\cdots<\varphi\left(  u\right)  $. Then, $\varphi
=\operatorname*{id}$.
\end{lemma}

\begin{proof}
[Proof of Lemma \ref{lem.sol.addexe.jacobi-complement.strinc}.]We have
$\varphi\left(  1\right)  <\varphi\left(  2\right)  <\cdots<\varphi\left(
u\right)  $. In other words, if $g$ and $h$ are two elements of $\left\{
1,2,\ldots,u\right\}  $ such that $g<h$, then%
\begin{equation}
\varphi\left(  g\right)  <\varphi\left(  h\right)  .
\label{pf.lem.sol.addexe.jacobi-complement.strinc.1}%
\end{equation}

Now, if $g$ and $h$ are two elements of $\left\{  1,2,\ldots,u\right\}  $
satisfying $\varphi\left(  g\right)  =\varphi\left(  h\right)  $, then
$g=h$\ \ \ \ \footnote{\textit{Proof.} Let $g$ and $h$ be two elements of
$\left\{  1,2,\ldots,u\right\}  $ satisfying $\varphi\left(  g\right)
=\varphi\left(  h\right)  $. We must prove that $g=h$.
\par
Indeed, assume the contrary (for the sake of contradiction). Thus, $g\neq h$.
Assume WLOG that $g\leq h$ (otherwise, we can simply swap $g$ with $h$). Then,
$g<h$ (since $g\leq h$ and $g\neq h$). Thus, $\varphi\left(  g\right)
<\varphi\left(  h\right)  $ (by
(\ref{pf.lem.sol.addexe.jacobi-complement.strinc.1})). This contradicts
$\varphi\left(  g\right)  =\varphi\left(  h\right)  $. This contradiction
proves that our assumption was wrong. Hence, $g=h$ is proven. Qed.}. In other
words, the map $\varphi$ is injective.

But $\left\{  1,2,\ldots,u\right\}  $ is a finite set and satisfies
$\left\vert \left\{  1,2,\ldots,u\right\}  \right\vert \geq\left\vert \left\{
1,2,\ldots,u\right\}  \right\vert $. Hence, Lemma
\ref{lem.jectivity.pigeon-inj} (applied to $\left\{  1,2,\ldots,u\right\}  $,
$\left\{  1,2,\ldots,u\right\}  $ and $\varphi$ instead of $U$, $V$ and $f$)
shows that we have the following logical equivalence:%
\[
\left(  \varphi\text{ is injective}\right)  \ \Longleftrightarrow\ \left(
\varphi\text{ is bijective}\right)  .
\]
Thus, $\varphi$ is bijective (since $\varphi$ is injective). Hence, $\varphi$
is a bijective map $\left\{  1,2,\ldots,u\right\}  \rightarrow\left\{
1,2,\ldots,u\right\}  $. In other words, $\varphi$ is a permutation of
$\left\{  1,2,\ldots,u\right\}  $. In other words, $\varphi\in S_{u}$ (because
$S_{u}$ is the set of all permutations of $\left\{  1,2,\ldots,u\right\}  $
(by the definition of $S_{u}$)).

Now, recall that $\varphi\left(  1\right)  <\varphi\left(  2\right)
<\cdots<\varphi\left(  u\right)  $. In other words, $\varphi\left(  k\right)
<\varphi\left(  k+1\right)  $ for each $k\in\left\{  1,2,\ldots,u-1\right\}
$. Hence, $\varphi\left(  k\right)  \leq\varphi\left(  k+1\right)  $ for each
$k\in\left\{  1,2,\ldots,u-1\right\}  $. In other words, $\varphi\left(
1\right)  \leq\varphi\left(  2\right)  \leq\cdots\leq\varphi\left(  u\right)
$. Hence, Proposition \ref{prop.sol.ps2.2.5.d} (applied to $n=u$ and
$\sigma=\varphi$) yields that $\varphi=\operatorname*{id}$. This proves Lemma
\ref{lem.sol.addexe.jacobi-complement.strinc}.
\end{proof}

\begin{lemma}
\label{lem.sol.addexe.jacobi-complement.sets-eq}Let $u\in\mathbb{N}$. Let
$\left(  p_{1},p_{2},\ldots,p_{u}\right)  $ be a list of integers such that
$p_{1}<p_{2}<\cdots<p_{u}$. Let $\left(  q_{1},q_{2},\ldots,q_{u}\right)  $ be
a list of integers such that $q_{1}<q_{2}<\cdots<q_{u}$. Assume that $\left\{
p_{1},p_{2},\ldots,p_{u}\right\}  =\left\{  q_{1},q_{2},\ldots,q_{u}\right\}
$. Then, $\left(  p_{1},p_{2},\ldots,p_{u}\right)  =\left(  q_{1},q_{2}%
,\ldots,q_{u}\right)  $.
\end{lemma}

\begin{proof}
[First proof of Lemma \ref{lem.sol.addexe.jacobi-complement.sets-eq}.]We have
$\left\{  p_{1},p_{2},\ldots,p_{u}\right\}  =\left\{  q_{1},q_{2},\ldots
,q_{u}\right\}  $. Thus, we can define a set $S$ by setting%
\[
S=\left\{  p_{1},p_{2},\ldots,p_{u}\right\}  =\left\{  q_{1},q_{2}%
,\ldots,q_{u}\right\}  .
\]
Consider this $S$. This set $S$ is finite (since $S=\left\{  p_{1}%
,p_{2},\ldots,p_{u}\right\}  $). Furthermore, $S=\left\{  p_{1},p_{2}%
,\ldots,p_{u}\right\}  $; thus, $S$ is a set of integers (since $p_{1}%
,p_{2},\ldots,p_{u}$ are integers\footnote{since $\left(  p_{1},p_{2}%
,\ldots,p_{u}\right)  $ is a list of integers}).

The elements $p_{1},p_{2},\ldots,p_{u}$ belong to $S$ (since $S=\left\{
p_{1},p_{2},\ldots,p_{u}\right\}  $). Hence, $\left(  p_{1},p_{2},\ldots
,p_{u}\right)  $ is a list of elements of $S$.

Now, recall Definition \ref{def.ind.inclist0}. The list $\left(  p_{1}%
,p_{2},\ldots,p_{u}\right)  $ is a list of elements of $S$ such that
$S=\left\{  p_{1},p_{2},\ldots,p_{u}\right\}  $ and $p_{1}<p_{2}<\cdots<p_{u}%
$. In other words, this list $\left(  p_{1},p_{2},\ldots,p_{u}\right)  $ is an
increasing list of $S$ (by the definition of an \textquotedblleft increasing
list\textquotedblright). The same argument (applied to $q_{1},q_{2}%
,\ldots,q_{u}$ instead of $p_{1},p_{2},\ldots,p_{u}$) shows that the list
$\left(  q_{1},q_{2},\ldots,q_{u}\right)  $ is an increasing list of $S$.

But Theorem \ref{thm.ind.inclist.unex} shows that $S$ has exactly one
increasing list. Hence, $S$ has \textbf{at most one} increasing list. In other
words, any two increasing lists of $S$ are equal. Thus, $\left(  p_{1}%
,p_{2},\ldots,p_{u}\right)  $ and $\left(  q_{1},q_{2},\ldots,q_{u}\right)  $
are equal (since $\left(  p_{1},p_{2},\ldots,p_{u}\right)  $ and $\left(
q_{1},q_{2},\ldots,q_{u}\right)  $ are two increasing lists of $S$). In other
words, $\left(  p_{1},p_{2},\ldots,p_{u}\right)  =\left(  q_{1},q_{2}%
,\ldots,q_{u}\right)  $. This proves Lemma
\ref{lem.sol.addexe.jacobi-complement.sets-eq}.
\end{proof}

\begin{proof}
[Second proof of Lemma \ref{lem.sol.addexe.jacobi-complement.sets-eq}.]We have
$p_{1}<p_{2}<\cdots<p_{u}$. In other words,%
\begin{equation}
p_{k}<p_{k+1}\ \ \ \ \ \ \ \ \ \ \text{for every }k\in\left\{  1,2,\ldots
,u-1\right\}  . \label{pf.lem.sol.addexe.jacobi-complement.sets-eq.p}%
\end{equation}

We have $q_{1}<q_{2}<\cdots<q_{u}$. Hence, if $g$ and $h$ are two elements of
$\left\{  1,2,\ldots,u\right\}  $ satisfying $g\leq h$, then%
\begin{equation}
q_{g}\leq q_{h}. \label{pf.lem.sol.addexe.jacobi-complement.sets-eq.q}%
\end{equation}

Define a map $\varphi:\left\{  1,2,\ldots,u\right\}  \rightarrow\left\{
1,2,\ldots,u\right\}  $ as follows:

Let $k\in\left\{  1,2,\ldots,u\right\}  $. Then, $p_{k}\in\left\{  p_{1}%
,p_{2},\ldots,p_{u}\right\}  =\left\{  q_{1},q_{2},\ldots,q_{u}\right\}  $.
Hence, there exists an $\ell\in\left\{  1,2,\ldots,u\right\}  $ such that
$p_{k}=q_{\ell}$. Consider this $\ell$. Define $\varphi\left(  k\right)  $ by
$\varphi\left(  k\right)  =\ell$. Note that
\begin{equation}
p_{k}=q_{\ell}=q_{\varphi\left(  k\right)  }
\label{pf.lem.sol.addexe.jacobi-complement.sets-eq.1}%
\end{equation}
(since $\ell=\varphi\left(  k\right)  $).

Now, forget that we fixed $k$. Thus, for each $k\in\left\{  1,2,\ldots
,u\right\}  $, we have defined an element $\varphi\left(  k\right)
\in\left\{  1,2,\ldots,u\right\}  $. In other words, we have defined a map
$\varphi:\left\{  1,2,\ldots,u\right\}  \rightarrow\left\{  1,2,\ldots
,u\right\}  $. This map furthermore satisfies
(\ref{pf.lem.sol.addexe.jacobi-complement.sets-eq.1}) for each $k\in\left\{
1,2,\ldots,u\right\}  $.

Now, every $k\in\left\{  1,2,\ldots,u-1\right\}  $ satisfies $\varphi\left(
k\right)  <\varphi\left(  k+1\right)  $\ \ \ \ \footnote{\textit{Proof.} Let
$k\in\left\{  1,2,\ldots,u-1\right\}  $. Then, $k\in\left\{  1,2,\ldots
,u-1\right\}  \subseteq\left\{  1,2,\ldots,u\right\}  $; thus, $\varphi\left(
k\right)  $ is well-defined.
\par
Also, from $k\in\left\{  1,2,\ldots,u-1\right\}  $, we obtain $k+1\in\left\{
2,3,\ldots,u\right\}  \subseteq\left\{  1,2,\ldots,u\right\}  $; thus,
$\varphi\left(  k+1\right)  $ is well-defined.
\par
From (\ref{pf.lem.sol.addexe.jacobi-complement.sets-eq.p}), we obtain
$p_{k}<p_{k+1}$. But (\ref{pf.lem.sol.addexe.jacobi-complement.sets-eq.1})
yields $p_{k}=q_{\varphi\left(  k\right)  }$ (since $k\in\left\{
1,2,\ldots,u\right\}  $). Thus, $q_{\varphi\left(  k\right)  }=p_{k}$. Also,
(\ref{pf.lem.sol.addexe.jacobi-complement.sets-eq.1}) (applied to $k+1$
instead of $k$) yields $p_{k+1}=q_{\varphi\left(  k+1\right)  }$ (since
$k+1\in\left\{  1,2,\ldots,u\right\}  $). Hence, $q_{\varphi\left(  k\right)
}=p_{k}<p_{k+1}=q_{\varphi\left(  k+1\right)  }$.
\par
Now, assume (for the sake of contradiction) that $\varphi\left(  k\right)
\geq\varphi\left(  k+1\right)  $. Thus, $\varphi\left(  k+1\right)
\leq\varphi\left(  k\right)  $. Hence,
(\ref{pf.lem.sol.addexe.jacobi-complement.sets-eq.q}) (applied to
$g=\varphi\left(  k+1\right)  $ and $h=\varphi\left(  k\right)  $) yields
$q_{\varphi\left(  k+1\right)  }\leq q_{\varphi\left(  k\right)  }%
<q_{\varphi\left(  k+1\right)  }$. This is absurd. Thus, we have found a
contradiction. Hence, our assumption (that $\varphi\left(  k\right)
\geq\varphi\left(  k+1\right)  $) must have been false. Therefore, we cannot
have $\varphi\left(  k\right)  \geq\varphi\left(  k+1\right)  $. Hence, we
must have $\varphi\left(  k\right)  <\varphi\left(  k+1\right)  $. Qed.}. In
other words, we have $\varphi\left(  1\right)  <\varphi\left(  2\right)
<\cdots<\varphi\left(  u\right)  $. Hence, Lemma
\ref{lem.sol.addexe.jacobi-complement.strinc} shows that $\varphi
=\operatorname*{id}$. Now, every $k\in\left\{  1,2,\ldots,u\right\}  $
satisfies%
\begin{align*}
p_{k}  &  =q_{\varphi\left(  k\right)  }\ \ \ \ \ \ \ \ \ \ \left(  \text{by
(\ref{pf.lem.sol.addexe.jacobi-complement.sets-eq.1})}\right) \\
&  =q_{k}\ \ \ \ \ \ \ \ \ \ \left(  \text{since }\underbrace{\varphi
}_{=\operatorname*{id}}\left(  k\right)  =\operatorname*{id}\left(  k\right)
=k\right)  .
\end{align*}
In other words, $\left(  p_{1},p_{2},\ldots,p_{u}\right)  =\left(  q_{1}%
,q_{2},\ldots,q_{u}\right)  $. This proves Lemma
\ref{lem.sol.addexe.jacobi-complement.sets-eq}.
\end{proof}

\begin{proof}
[Proof of Lemma \ref{lem.sol.addexe.jacobi-complement.increase-bij}.]We have
$S\subseteq\mathbb{Z}$ (since $S$ is a set of integers).

The definition of $\mathbf{I}$ yields%
\[
\mathbf{I}=\left\{  \left(  g_{1},g_{2},\ldots,g_{u}\right)  \in S^{u}%
\ \mid\ g_{1}<g_{2}<\cdots<g_{u}\right\}  .
\]
Hence, every element of $\mathbf{I}$ has the form $\left(  g_{1},g_{2}%
,\ldots,g_{u}\right)  $ for some $\left(  g_{1},g_{2},\ldots,g_{u}\right)  \in
S^{u}$.

For every $R\in\mathcal{P}_{u}\left(  S\right)  $, the list $w\left(
R\right)  $ is well-defined\footnote{\textit{Proof.} Let $R\in\mathcal{P}%
_{u}\left(  S\right)  $. We must show that the list $w\left(  R\right)  $ is
well-defined.
\par
We have $R\in\mathcal{P}_{u}\left(  S\right)  $. In other words, $R$ is a
$u$-element subset of $S$ (since $\mathcal{P}_{u}\left(  S\right)  $ is the
set of all $u$-element subsets of $S$). Now, $R$ is a subset of $S$; in other
words, $R\subseteq S$. Hence, $R\subseteq S\subseteq\mathbb{Z}$ (since $S$ is
a set of integers). In other words, $R$ is a set of integers. Also, $R$ is
finite (since $R$ is a $u$-element set). Hence, $R$ is a finite set of
integers.
\par
But the list $w\left(  I\right)  $ is well-defined whenever $I$ is a finite
set of integers. Applying this to $I=R$, we conclude that the list $w\left(
R\right)  $ is well-defined. Qed.} and belongs to $\mathbf{I}$%
\ \ \ \ \footnote{\textit{Proof.} Let $R\in\mathcal{P}_{u}\left(  S\right)  $.
We must show that the list $w\left(  R\right)  $ belongs to $\mathbf{I}$.
\par
We have $R\in\mathcal{P}_{u}\left(  S\right)  $. In other words, $R$ is a
$u$-element subset of $S$ (since $\mathcal{P}_{u}\left(  S\right)  $ is the
set of all $u$-element subsets of $S$). Now, $R$ is a subset of $S$; in other
words, $R\subseteq S$. Hence, $R\subseteq S\subseteq\mathbb{Z}$ (since $S$ is
a set of integers). In other words, $R$ is a set of integers. Also, $R$ is
finite (since $R$ is a $u$-element set). Hence, $R$ is a finite set of
integers. Also, $\left\vert R\right\vert =u$ (since $R$ is a $u$-element set).
\par
The definition of $w\left(  R\right)  $ shows that $w\left(  R\right)  $ is
the list of all elements of $R$ in increasing order (with no repetitions).
Hence, this list $w\left(  R\right)  $ has $\left\vert R\right\vert $ entries.
In other words, the list $w\left(  R\right)  $ has $u$ entries (since
$\left\vert R\right\vert =u$).
\par
Write the list $w\left(  R\right)  $ in the form $w\left(  R\right)  =\left(
r_{1},r_{2},\ldots,r_{u}\right)  $. (This is possible, since the list
$w\left(  R\right)  $ has $u$ entries.)
\par
The list $w\left(  R\right)  $ is the list of all elements of $R$ in
increasing order (with no repetitions). Hence, $w\left(  R\right)  $ is a list
of all elements of $R$. In other words, $\left(  r_{1},r_{2},\ldots
,r_{u}\right)  $ is a list of all elements of $R$ (since $w\left(  R\right)
=\left(  r_{1},r_{2},\ldots,r_{u}\right)  $). Hence, $\left\{  r_{1}%
,r_{2},\ldots,r_{u}\right\}  =R$. Now, each $i\in\left\{  1,2,\ldots
,u\right\}  $ satisfies $r_{i}\in\left\{  r_{1},r_{2},\ldots,r_{u}\right\}
=R\subseteq S$. Thus, $\left(  r_{1},r_{2},\ldots,r_{u}\right)  \in S^{u}$.
\par
Moreover, the list $w\left(  R\right)  $ is the list of all elements of $R$ in
increasing order (with no repetitions). Hence, the list $w\left(  R\right)  $
is strictly increasing. In other words, the list $\left(  r_{1},r_{2}%
,\ldots,r_{u}\right)  $ is strictly increasing (since $w\left(  R\right)
=\left(  r_{1},r_{2},\ldots,r_{u}\right)  $). In other words, $r_{1}%
<r_{2}<\cdots<r_{u}$.
\par
Now, we know that the list $\left(  r_{1},r_{2},\ldots,r_{u}\right)  \in
S^{u}$ satisfies $r_{1}<r_{2}<\cdots<r_{u}$ and $w\left(  R\right)  =\left(
r_{1},r_{2},\ldots,r_{u}\right)  $. Hence, there exists some $\left(
g_{1},g_{2},\ldots,g_{u}\right)  \in S^{u}$ satisfying $g_{1}<g_{2}%
<\cdots<g_{u}$ and $w\left(  R\right)  =\left(  g_{1},g_{2},\ldots
,g_{u}\right)  $ (namely, $\left(  g_{1},g_{2},\ldots,g_{u}\right)  =\left(
r_{1},r_{2},\ldots,r_{u}\right)  $). In other words,%
\[
w\left(  R\right)  \in\left\{  \left(  g_{1},g_{2},\ldots,g_{u}\right)  \in
S^{u}\ \mid\ g_{1}<g_{2}<\cdots<g_{u}\right\}  .
\]
This rewrites as $w\left(  R\right)  \in\mathbf{I}$ (since $\mathbf{I}%
=\left\{  \left(  g_{1},g_{2},\ldots,g_{u}\right)  \in S^{u}\ \mid
\ g_{1}<g_{2}<\cdots<g_{u}\right\}  $). In other words, $w\left(  R\right)  $
belongs to $\mathbf{I}$. Qed.}. Hence, the map%
\begin{align*}
\mathcal{P}_{u}\left(  S\right)   &  \rightarrow\mathbf{I},\\
R  &  \mapsto w\left(  R\right)
\end{align*}
is well-defined. Denote this map by $\alpha$.

Now, every $\left(  r_{1},r_{2},\ldots,r_{u}\right)  \in\mathbf{I}$ satisfies
$\left\{  r_{1},r_{2},\ldots,r_{u}\right\}  \in\mathcal{P}_{u}\left(
S\right)  $\ \ \ \ \footnote{\textit{Proof.} Let $\left(  r_{1},r_{2}%
,\ldots,r_{u}\right)  \in\mathbf{I}$. Then,%
\[
\left(  r_{1},r_{2},\ldots,r_{u}\right)  \in\mathbf{I}=\left\{  \left(
g_{1},g_{2},\ldots,g_{u}\right)  \in S^{u}\ \mid\ g_{1}<g_{2}<\cdots
<g_{u}\right\}  .
\]
In other words, $\left(  r_{1},r_{2},\ldots,r_{u}\right)  $ is an element
$\left(  g_{1},g_{2},\ldots,g_{u}\right)  \in S^{u}$ satisfying $g_{1}%
<g_{2}<\cdots<g_{u}$. In other words, $\left(  r_{1},r_{2},\ldots
,r_{u}\right)  $ is an element of $S^{u}$ and satisfies $r_{1}<r_{2}%
<\cdots<r_{u}$.
\par
We have $\left(  r_{1},r_{2},\ldots,r_{u}\right)  \in S^{u}$. In other words,
$r_{i}\in S$ for each $i\in\left\{  1,2,\ldots,u\right\}  $. Hence, $\left\{
r_{1},r_{2},\ldots,r_{u}\right\}  \subseteq S$. In other words, $\left\{
r_{1},r_{2},\ldots,r_{u}\right\}  $ is a subset of $S$.
\par
The $u$ elements $r_{1},r_{2},\ldots,r_{u}$ are pairwise distinct (since
$r_{1}<r_{2}<\cdots<r_{u}$). Hence, $\left\vert \left\{  r_{1},r_{2}%
,\ldots,r_{u}\right\}  \right\vert =u$. Thus, $\left\{  r_{1},r_{2}%
,\ldots,r_{u}\right\}  $ is a $u$-element set. Now, we know that $\left\{
r_{1},r_{2},\ldots,r_{u}\right\}  $ is a $u$-element subset of $S$ (since
$\left\{  r_{1},r_{2},\ldots,r_{u}\right\}  $ is a $u$-element set and a
subset of $S$). In other words, $\left\{  r_{1},r_{2},\ldots,r_{u}\right\}
\in\mathcal{P}_{u}\left(  S\right)  $ (since $\mathcal{P}_{u}\left(  S\right)
$ is the set of all $u$-element subsets of $S$). Qed.}.

Recall that every element of $\mathbf{I}$ has the form $\left(  g_{1}%
,g_{2},\ldots,g_{u}\right)  $ for some $\left(  g_{1},g_{2},\ldots
,g_{u}\right)  \in S^{u}$. In other words, every element of $\mathbf{I}$ has
the form $\left(  r_{1},r_{2},\ldots,r_{u}\right)  $ for some $\left(
r_{1},r_{2},\ldots,r_{u}\right)  \in S^{u}$ (here, we have renamed the index
$\left(  g_{1},g_{2},\ldots,g_{u}\right)  $ as $\left(  r_{1},r_{2}%
,\ldots,r_{u}\right)  $). Hence, we can define a map $\beta:\mathbf{I}%
\rightarrow\mathcal{P}_{u}\left(  S\right)  $ by%
\[
\left(  \beta\left(  r_{1},r_{2},\ldots,r_{u}\right)  =\left\{  r_{1}%
,r_{2},\ldots,r_{u}\right\}  \ \ \ \ \ \ \ \ \ \ \text{for every }\left(
r_{1},r_{2},\ldots,r_{u}\right)  \in\mathbf{I}\right)
\]
(since every $\left(  r_{1},r_{2},\ldots,r_{u}\right)  \in\mathbf{I}$
satisfies $\left\{  r_{1},r_{2},\ldots,r_{u}\right\}  \in\mathcal{P}%
_{u}\left(  S\right)  $). Consider this $\beta$.

We have $\alpha\circ\beta=\operatorname*{id}$\ \ \ \ \footnote{\textit{Proof.}
Let $\mathbf{i}\in\mathbf{I}$. Then,%
\[
\mathbf{i}\in\mathbf{I}=\left\{  \left(  g_{1},g_{2},\ldots,g_{u}\right)  \in
S^{u}\ \mid\ g_{1}<g_{2}<\cdots<g_{u}\right\}  .
\]
In other words, $\mathbf{i}$ has the form $\mathbf{i}=\left(  g_{1}%
,g_{2},\ldots,g_{u}\right)  $ for some $\left(  g_{1},g_{2},\ldots
,g_{u}\right)  \in S^{u}$ satisfying $g_{1}<g_{2}<\cdots<g_{u}$. Consider this
$\left(  g_{1},g_{2},\ldots,g_{u}\right)  $.
\par
Set $I=\beta\left(  \mathbf{i}\right)  $. Notice that $\beta\left(
\mathbf{i}\right)  \in\mathcal{P}_{u}\left(  S\right)  $ (because $\beta$ is a
map $\mathbf{I}\rightarrow\mathcal{P}_{u}\left(  S\right)  $). In other words,
$I\in\mathcal{P}_{u}\left(  S\right)  $ (since $I=\beta\left(  \mathbf{i}%
\right)  $). In other words, $I$ is a $u$-element subset of $S$ (since
$\mathcal{P}_{u}\left(  S\right)  $ is the set of all $u$-element subsets of
$S$). Thus, $I$ is a $u$-element set. In other words, $\left\vert I\right\vert
=u$.
\par
Applying the map $\beta$ to the equality $\mathbf{i}=\left(  g_{1}%
,g_{2},\ldots,g_{u}\right)  $, we obtain%
\[
\beta\left(  \mathbf{i}\right)  =\beta\left(  g_{1},g_{2},\ldots,g_{u}\right)
=\left\{  g_{1},g_{2},\ldots,g_{u}\right\}
\]
(by the definition of $\beta$). Thus, $I=\beta\left(  \mathbf{i}\right)
=\left\{  g_{1},g_{2},\ldots,g_{u}\right\}  $.
\par
Now, $\alpha\left(  \underbrace{\beta\left(  \mathbf{i}\right)  }_{=I}\right)
=\alpha\left(  I\right)  =w\left(  I\right)  $ (by the definition of the map
$\alpha$). But $w\left(  I\right)  $ is the list of all elements of $I$ in
increasing order (with no repetitions) (by the definition of $w\left(
I\right)  $). Hence, $w\left(  I\right)  $ is a list of size $\left\vert
I\right\vert =u$.
\par
Write the list $w\left(  I\right)  $ in the form $w\left(  I\right)  =\left(
p_{1},p_{2},\ldots,p_{u}\right)  $. (This is possible, since $w\left(
I\right)  $ is a list of size $u$.)
\par
Recall that $w\left(  I\right)  $ is the list of all elements of $I$ in
increasing order (with no repetitions). In other words, $\left(  p_{1}%
,p_{2},\ldots,p_{u}\right)  $ is the list of all elements of $I$ in increasing
order (with no repetitions) (since $w\left(  I\right)  =\left(  p_{1}%
,p_{2},\ldots,p_{u}\right)  $). Thus, $\left(  p_{1},p_{2},\ldots
,p_{u}\right)  $ is a strictly increasing list. In other words, $p_{1}%
<p_{2}<\cdots<p_{u}$.
\par
But $\left(  p_{1},p_{2},\ldots,p_{u}\right)  $ is the list of all elements of
$I$ in increasing order (with no repetitions). Thus, $\left(  p_{1}%
,p_{2},\ldots,p_{u}\right)  $ is a list of all elements of $I$. Hence,
$\left\{  p_{1},p_{2},\ldots,p_{u}\right\}  =I$. Thus, $\left\{  p_{1}%
,p_{2},\ldots,p_{u}\right\}  =I=\left\{  g_{1},g_{2},\ldots,g_{u}\right\}  $.
\par
Note that $I$ is a subset of $S$. Thus, $I\subseteq S\subseteq\mathbb{Z}$.
Each $k\in\left\{  1,2,\ldots,u\right\}  $ satisfies $p_{k}\in\left\{
p_{1},p_{2},\ldots,p_{u}\right\}  =I\subseteq\mathbb{Z}$. In other words, for
each $k\in\left\{  1,2,\ldots,u\right\}  $, the element $p_{k}$ is an integer.
Hence, $\left(  p_{1},p_{2},\ldots,p_{u}\right)  $ is a list of integers.
Also, each $k\in\left\{  1,2,\ldots,u\right\}  $ satisfies $g_{k}\in\left\{
g_{1},g_{2},\ldots,g_{u}\right\}  =I\subseteq\mathbb{Z}$. In other words, for
each $k\in\left\{  1,2,\ldots,u\right\}  $, the element $g_{k}$ is an integer.
Hence, $\left(  g_{1},g_{2},\ldots,g_{u}\right)  $ is a list of integers.
\par
Now, Lemma \ref{lem.sol.addexe.jacobi-complement.sets-eq} (applied to $\left(
g_{1},g_{2},\ldots,g_{u}\right)  $ instead of $\left(  q_{1},q_{2}%
,\ldots,q_{u}\right)  $) yields that $\left(  p_{1},p_{2},\ldots,p_{u}\right)
=\left(  g_{1},g_{2},\ldots,g_{u}\right)  $.
\par
But%
\begin{align*}
\left(  \alpha\circ\beta\right)  \left(  \mathbf{i}\right)   &  =\alpha\left(
\beta\left(  \mathbf{i}\right)  \right)  =w\left(  I\right)  =\left(
p_{1},p_{2},\ldots,p_{u}\right)  =\left(  g_{1},g_{2},\ldots,g_{u}\right) \\
&  =\mathbf{i}\ \ \ \ \ \ \ \ \ \ \left(  \text{since }\mathbf{i}=\left(
g_{1},g_{2},\ldots,g_{u}\right)  \right) \\
&  =\operatorname*{id}\left(  \mathbf{i}\right)  .
\end{align*}
\par
Now, forget that we fixed $\mathbf{i}$. We thus have shown that $\left(
\alpha\circ\beta\right)  \left(  \mathbf{i}\right)  =\operatorname*{id}\left(
\mathbf{i}\right)  $ for each $\mathbf{i}\in\mathbf{I}$. In other words,
$\alpha\circ\beta=\operatorname*{id}$. Qed.} and $\beta\circ\alpha
=\operatorname*{id}$\ \ \ \ \footnote{\textit{Proof.} Let $I\in\mathcal{P}%
_{u}\left(  S\right)  $. Then, $\alpha\left(  I\right)  =w\left(  I\right)  $
(by the definition of $\alpha$). Clearly, $\alpha\left(  I\right)
\in\mathbf{I}$ (since $\alpha$ is a map $\mathcal{P}_{u}\left(  S\right)
\rightarrow\mathbf{I}$).
\par
We have $I\in\mathcal{P}_{u}\left(  S\right)  $. In other words, $I$ is a
$u$-element subset of $S$ (since $\mathcal{P}_{u}\left(  S\right)  $ is the
set of all $u$-element subsets of $S$). Thus, $I$ is a $u$-element set. In
other words, $\left\vert I\right\vert =u$.
\par
But $w\left(  I\right)  $ is the list of all elements of $I$ in increasing
order (with no repetitions) (by the definition of $w\left(  I\right)  $).
Hence, $w\left(  I\right)  $ is a list of size $\left\vert I\right\vert =u$.
In other words, $\alpha\left(  I\right)  $ is a list of size $u$ (since
$\alpha\left(  I\right)  =w\left(  I\right)  $).
\par
Write the list $\alpha\left(  I\right)  $ in the form $\alpha\left(  I\right)
=\left(  p_{1},p_{2},\ldots,p_{u}\right)  $. (This is possible, since
$\alpha\left(  I\right)  $ is a list of size $u$.) Hence, $w\left(  I\right)
=\alpha\left(  I\right)  =\left(  p_{1},p_{2},\ldots,p_{u}\right)  $.
\par
Recall that $w\left(  I\right)  $ is the list of all elements of $I$ in
increasing order (with no repetitions). Thus, $w\left(  I\right)  $ is a list
of all elements of $I$. In other words, $\left(  p_{1},p_{2},\ldots
,p_{u}\right)  $ is a list of all elements of $I$ (since $w\left(  I\right)
=\left(  p_{1},p_{2},\ldots,p_{u}\right)  $). Thus, $\left\{  p_{1}%
,p_{2},\ldots,p_{u}\right\}  =I$.
\par
But%
\begin{align*}
\left(  \beta\circ\alpha\right)  \left(  I\right)   &  =\beta\left(
\underbrace{\alpha\left(  I\right)  }_{=\left(  p_{1},p_{2},\ldots
,p_{u}\right)  }\right)  =\beta\left(  p_{1},p_{2},\ldots,p_{u}\right) \\
&  =\left\{  p_{1},p_{2},\ldots,p_{u}\right\}  \ \ \ \ \ \ \ \ \ \ \left(
\text{by the definition of }\beta\right) \\
&  =I=\operatorname*{id}\left(  I\right)  .
\end{align*}
\par
Let us now forget that we fixed $I$. We thus have shown that $\left(
\beta\circ\alpha\right)  \left(  I\right)  =\operatorname*{id}\left(
I\right)  $ for each $I\in\mathcal{P}_{u}\left(  S\right)  $. In other words,
$\beta\circ\alpha=\operatorname*{id}$. Qed.}. Combining these two equalities,
we conclude that the maps $\alpha$ and $\beta$ are mutually inverse. Hence,
the map $\alpha$ is invertible. In other words, the map $\alpha$ is a
bijection. In other words, the map%
\begin{align*}
\mathcal{P}_{u}\left(  S\right)   &  \rightarrow\mathbf{I},\\
R  &  \mapsto w\left(  R\right)
\end{align*}
is a bijection (since the map $\alpha$ is the map%
\begin{align*}
\mathcal{P}_{u}\left(  S\right)   &  \rightarrow\mathbf{I},\\
R  &  \mapsto w\left(  R\right)
\end{align*}
(by the definition of $\alpha$)). This completes the proof of Lemma
\ref{lem.sol.addexe.jacobi-complement.increase-bij}.
\end{proof}
\end{verlong}

We can use Lemma \ref{lem.sol.addexe.jacobi-complement.increase-bij} to obtain
the following (slightly less general) form of Corollary
\ref{cor.adj(AB).cauchy-binet-general}:

\begin{corollary}
\label{cor.sol.addexe.jacobi-complement.CB}Let $n\in\mathbb{N}$,
$m\in\mathbb{N}$ and $p\in\mathbb{N}$. Let $A$ be an $n\times p$-matrix. Let
$B$ be a $p\times m$-matrix. Let $k\in\mathbb{N}$. Let $P$ be a subset of
$\left\{  1,2,\ldots,n\right\}  $ such that $\left\vert P\right\vert =k$. Let
$Q$ be a subset of $\left\{  1,2,\ldots,m\right\}  $ such that $\left\vert
Q\right\vert =k$. Then,%
\[
\det\left(  \operatorname*{sub}\nolimits_{w\left(  P\right)  }^{w\left(
Q\right)  }\left(  AB\right)  \right)  =\sum_{\substack{R\subseteq\left\{
1,2,\ldots,p\right\}  ;\\\left\vert R\right\vert =k}}\det\left(
\operatorname*{sub}\nolimits_{w\left(  P\right)  }^{w\left(  R\right)
}A\right)  \cdot\det\left(  \operatorname*{sub}\nolimits_{w\left(  R\right)
}^{w\left(  Q\right)  }B\right)  .
\]

\end{corollary}

\begin{vershort}
\begin{proof}
[Proof of Corollary \ref{cor.sol.addexe.jacobi-complement.CB}.]For every set
$S$, we let $\mathcal{P}_{k}\left(  S\right)  $ denote the set of all
$k$-element subsets of $S$.

Recall that $w\left(  P\right)  $ was defined as the list of all elements of
$P$ in increasing order (with no repetitions). Thus, $w\left(  P\right)  $ is
a list of $\left\vert P\right\vert $ elements. In other words, $w\left(
P\right)  $ is a list of $k$ elements (since $\left\vert P\right\vert =k$).
Similarly, $w\left(  Q\right)  $ is a list of $k$ elements.

Write the list $w\left(  P\right)  $ in the form $w\left(  P\right)  =\left(
i_{1},i_{2},\ldots,i_{k}\right)  $. (This is possible, since $w\left(
P\right)  $ is a list of $k$ elements.)

Write the list $w\left(  Q\right)  $ in the form $w\left(  Q\right)  =\left(
j_{1},j_{2},\ldots,j_{k}\right)  $. (This is possible, since $w\left(
Q\right)  $ is a list of $k$ elements.)

We have $\left(  i_{1},i_{2},\ldots,i_{k}\right)  =w\left(  P\right)  $. Thus,
the elements $i_{1},i_{2},\ldots,i_{k}$ are elements of $P$, and hence also
elements of $\left\{  1,2,\ldots,n\right\}  $ (since $P\subseteq\left\{
1,2,\ldots,n\right\}  $). The same argument (but applied to $m$, $Q$ and
$\left(  j_{1},j_{2},\ldots,j_{k}\right)  $ instead of $n$, $P$ and $\left(
i_{1},i_{2},\ldots,i_{k}\right)  $) yields that $j_{1},j_{2},\ldots,j_{k}$ are
elements of $\left\{  1,2,\ldots,m\right\}  $. Hence, Corollary
\ref{cor.adj(AB).cauchy-binet-general} (applied to $u=k$) yields%
\begin{align}
&  \det\left(  \operatorname*{sub}\nolimits_{i_{1},i_{2},\ldots,i_{k}}%
^{j_{1},j_{2},\ldots,j_{k}}\left(  AB\right)  \right) \nonumber\\
&  =\sum_{1\leq g_{1}<g_{2}<\cdots<g_{k}\leq p}\det\left(  \operatorname*{sub}%
\nolimits_{i_{1},i_{2},\ldots,i_{k}}^{g_{1},g_{2},\ldots,g_{k}}A\right)
\cdot\det\left(  \operatorname*{sub}\nolimits_{g_{1},g_{2},\ldots,g_{k}%
}^{j_{1},j_{2},\ldots,j_{k}}B\right)  .
\label{pf.cor.sol.addexe.jacobi-complement.CB.short.1}%
\end{align}
(Here, the summation sign \textquotedblleft$\sum_{1\leq g_{1}<g_{2}%
<\cdots<g_{k}\leq p}$\textquotedblright\ has to be interpreted as
\textquotedblleft$\sum_{\substack{\left(  g_{1},g_{2},\ldots,g_{k}\right)
\in\left\{  1,2,\ldots,p\right\}  ^{k};\\g_{1}<g_{2}<\cdots<g_{k}}%
}$\textquotedblright, in analogy to Remark \ref{rmk.cauchy-binet.sumsign}.)

Now, let $S$ denote the set $\left\{  1,2,\ldots,p\right\}  $. Recall that
$\mathcal{P}_{k}\left(  S\right)  $ denotes the set of all $k$-element subsets
of $S$. In other words,
\[
\mathcal{P}_{k}\left(  S\right)  =\left\{  T\subseteq S\ \mid\ \left\vert
T\right\vert =k\right\}  .
\]

Let $\mathbf{I}$ denote the set
\[
\left\{  \left(  g_{1},g_{2},\ldots,g_{k}\right)  \in S^{k}\ \mid\ g_{1}%
<g_{2}<\cdots<g_{k}\right\}  .
\]
Lemma \ref{lem.sol.addexe.jacobi-complement.increase-bij} (applied to $u=k$)
shows that the map%
\begin{align*}
\mathcal{P}_{k}\left(  S\right)   &  \rightarrow\mathbf{I},\\
R  &  \mapsto w\left(  R\right)
\end{align*}
is well-defined and a bijection.

The definition of $\mathbf{I}$ yields $\mathbf{I}=\left\{  \left(  g_{1}%
,g_{2},\ldots,g_{k}\right)  \in S^{k}\ \mid\ g_{1}<g_{2}<\cdots<g_{k}\right\}
$. Hence, we have the following equality of summation signs:%
\[
\sum_{\left(  g_{1},g_{2},\ldots,g_{k}\right)  \in\mathbf{I}}=\sum
_{\substack{\left(  g_{1},g_{2},\ldots,g_{k}\right)  \in S^{k};\\g_{1}%
<g_{2}<\cdots<g_{k}}}=\sum_{\substack{\left(  g_{1},g_{2},\ldots,g_{k}\right)
\in\left\{  1,2,\ldots,p\right\}  ^{k};\\g_{1}<g_{2}<\cdots<g_{k}}}
\]
(since $S=\left\{  1,2,\ldots,p\right\}  $). Comparing this with%
\[
\sum_{1\leq g_{1}<g_{2}<\cdots<g_{k}\leq p}=\sum_{\substack{\left(
g_{1},g_{2},\ldots,g_{k}\right)  \in\left\{  1,2,\ldots,p\right\}
^{k};\\g_{1}<g_{2}<\cdots<g_{k}}},
\]
we obtain%
\[
\sum_{1\leq g_{1}<g_{2}<\cdots<g_{k}\leq p}=\sum_{\left(  g_{1},g_{2}%
,\ldots,g_{k}\right)  \in\mathbf{I}}.
\]
Now, (\ref{pf.cor.sol.addexe.jacobi-complement.CB.short.1}) becomes%
\begin{align*}
&  \det\left(  \operatorname*{sub}\nolimits_{i_{1},i_{2},\ldots,i_{k}}%
^{j_{1},j_{2},\ldots,j_{k}}\left(  AB\right)  \right) \\
&  =\underbrace{\sum_{1\leq g_{1}<g_{2}<\cdots<g_{k}\leq p}}_{=\sum_{\left(
g_{1},g_{2},\ldots,g_{k}\right)  \in\mathbf{I}}}\det\left(
\underbrace{\operatorname*{sub}\nolimits_{i_{1},i_{2},\ldots,i_{k}}%
^{g_{1},g_{2},\ldots,g_{k}}A}_{=\operatorname*{sub}\nolimits_{\left(
i_{1},i_{2},\ldots,i_{k}\right)  }^{\left(  g_{1},g_{2},\ldots,g_{k}\right)
}A}\right)  \cdot\det\left(  \underbrace{\operatorname*{sub}\nolimits_{g_{1}%
,g_{2},\ldots,g_{k}}^{j_{1},j_{2},\ldots,j_{k}}B}_{=\operatorname*{sub}%
\nolimits_{\left(  g_{1},g_{2},\ldots,g_{k}\right)  }^{\left(  j_{1}%
,j_{2},\ldots,j_{k}\right)  }B}\right) \\
&  =\sum_{\left(  g_{1},g_{2},\ldots,g_{k}\right)  \in\mathbf{I}}\det\left(
\operatorname*{sub}\nolimits_{\left(  i_{1},i_{2},\ldots,i_{k}\right)
}^{\left(  g_{1},g_{2},\ldots,g_{k}\right)  }A\right)  \cdot\det\left(
\operatorname*{sub}\nolimits_{\left(  g_{1},g_{2},\ldots,g_{k}\right)
}^{\left(  j_{1},j_{2},\ldots,j_{k}\right)  }B\right) \\
&  =\underbrace{\sum_{R\in\mathcal{P}_{k}\left(  S\right)  }}_{\substack{=\sum
_{\substack{R\subseteq S;\\\left\vert R\right\vert =k}}\\\text{(since
}\mathcal{P}_{k}\left(  S\right)  =\left\{  T\subseteq S\ \mid\ \left\vert
T\right\vert =k\right\}  \text{)}}}\det\left(  \underbrace{\operatorname*{sub}%
\nolimits_{\left(  i_{1},i_{2},\ldots,i_{k}\right)  }^{w\left(  R\right)  }%
A}_{\substack{=\operatorname*{sub}\nolimits_{w\left(  P\right)  }^{w\left(
R\right)  }A\\\text{(since }\left(  i_{1},i_{2},\ldots,i_{k}\right)  =w\left(
P\right)  \text{)}}}\right)  \cdot\det\left(  \underbrace{\operatorname*{sub}%
\nolimits_{w\left(  R\right)  }^{\left(  j_{1},j_{2},\ldots,j_{k}\right)  }%
B}_{\substack{=\operatorname*{sub}\nolimits_{w\left(  R\right)  }^{w\left(
Q\right)  }B\\\text{(since }\left(  j_{1},j_{2},\ldots,j_{k}\right)  =w\left(
Q\right)  \text{)}}}\right) \\
&  \ \ \ \ \ \ \ \ \ \ \ \ \ \ \ \ \ \ \ \ \left(
\begin{array}
[c]{c}%
\text{here, we have substituted }w\left(  R\right)  \text{ for }\left(
g_{1},g_{2},\ldots,g_{k}\right)  \text{, since}\\
\text{the map }\mathcal{P}_{k}\left(  S\right)  \rightarrow\mathbf{I}%
,\ R\mapsto w\left(  R\right)  \text{ is a bijection}%
\end{array}
\right) \\
&  =\underbrace{\sum_{\substack{R\subseteq S;\\\left\vert R\right\vert =k}%
}}_{\substack{=\sum_{\substack{R\subseteq\left\{  1,2,\ldots,p\right\}
;\\\left\vert R\right\vert =k}}\\\text{(since }S=\left\{  1,2,\ldots
,p\right\}  \text{)}}}\det\left(  \operatorname*{sub}\nolimits_{w\left(
P\right)  }^{w\left(  R\right)  }A\right)  \cdot\det\left(
\operatorname*{sub}\nolimits_{w\left(  R\right)  }^{w\left(  Q\right)
}B\right) \\
&  =\sum_{\substack{R\subseteq\left\{  1,2,\ldots,p\right\}  ;\\\left\vert
R\right\vert =k}}\det\left(  \operatorname*{sub}\nolimits_{w\left(  P\right)
}^{w\left(  R\right)  }A\right)  \cdot\det\left(  \operatorname*{sub}%
\nolimits_{w\left(  R\right)  }^{w\left(  Q\right)  }B\right)  .
\end{align*}
Comparing this with%
\begin{align*}
\det\left(  \underbrace{\operatorname*{sub}\nolimits_{i_{1},i_{2},\ldots
,i_{k}}^{j_{1},j_{2},\ldots,j_{k}}\left(  AB\right)  }_{=\operatorname*{sub}%
\nolimits_{\left(  i_{1},i_{2},\ldots,i_{k}\right)  }^{\left(  j_{1}%
,j_{2},\ldots,j_{k}\right)  }\left(  AB\right)  }\right)   &  =\det\left(
\operatorname*{sub}\nolimits_{\left(  i_{1},i_{2},\ldots,i_{k}\right)
}^{\left(  j_{1},j_{2},\ldots,j_{k}\right)  }\left(  AB\right)  \right)
=\det\left(  \operatorname*{sub}\nolimits_{w\left(  P\right)  }^{w\left(
Q\right)  }\left(  AB\right)  \right) \\
&  \ \ \ \ \ \ \ \ \ \ \left(
\begin{array}
[c]{c}%
\text{since }\left(  i_{1},i_{2},\ldots,i_{k}\right)  =w\left(  P\right) \\
\text{and }\left(  j_{1},j_{2},\ldots,j_{k}\right)  =w\left(  Q\right)
\end{array}
\right)  ,
\end{align*}
we obtain%
\[
\det\left(  \operatorname*{sub}\nolimits_{w\left(  P\right)  }^{w\left(
Q\right)  }\left(  AB\right)  \right)  =\sum_{\substack{R\subseteq\left\{
1,2,\ldots,p\right\}  ;\\\left\vert R\right\vert =k}}\det\left(
\operatorname*{sub}\nolimits_{w\left(  P\right)  }^{w\left(  R\right)
}A\right)  \cdot\det\left(  \operatorname*{sub}\nolimits_{w\left(  R\right)
}^{w\left(  Q\right)  }B\right)  .
\]
This proves Corollary \ref{cor.sol.addexe.jacobi-complement.CB}.
\end{proof}
\end{vershort}

\begin{verlong}
\begin{proof}
[Proof of Corollary \ref{cor.sol.addexe.jacobi-complement.CB}.]For every set
$S$, we let $\mathcal{P}_{k}\left(  S\right)  $ denote the set of all
$k$-element subsets of $S$.

We know that $w\left(  P\right)  $ is the list of all elements of $P$ in
increasing order (with no repetitions) (by the definition of $w\left(
P\right)  $). Thus, $w\left(  P\right)  $ is a list of $\left\vert
P\right\vert $ elements. In other words, $w\left(  P\right)  $ is a list of
$k$ elements (since $\left\vert P\right\vert =k$).

Write the list $w\left(  P\right)  $ in the form $w\left(  P\right)  =\left(
i_{1},i_{2},\ldots,i_{k}\right)  $. (This is possible, since $w\left(
P\right)  $ is a list of $k$ elements.)

We know that $w\left(  Q\right)  $ is the list of all elements of $Q$ in
increasing order (with no repetitions) (by the definition of $w\left(
Q\right)  $). Thus, $w\left(  Q\right)  $ is a list of $\left\vert
Q\right\vert $ elements. In other words, $w\left(  Q\right)  $ is a list of
$k$ elements (since $\left\vert Q\right\vert =k$).

Write the list $w\left(  Q\right)  $ in the form $w\left(  Q\right)  =\left(
j_{1},j_{2},\ldots,j_{k}\right)  $. (This is possible, since $w\left(
Q\right)  $ is a list of $k$ elements.)

We know that $w\left(  P\right)  $ is the list of all elements of $P$ in
increasing order (with no repetitions). Hence, $w\left(  P\right)  $ is a list
of all elements of $P$. In other words, $\left(  i_{1},i_{2},\ldots
,i_{k}\right)  $ is a list of all elements of $P$ (since $w\left(  P\right)
=\left(  i_{1},i_{2},\ldots,i_{k}\right)  $). Thus, $\left\{  i_{1}%
,i_{2},\ldots,i_{k}\right\}  =P$. Now, every $\ell\in\left\{  1,2,\ldots
,k\right\}  $ satisfies $i_{\ell}\in\left\{  i_{1},i_{2},\ldots,i_{k}\right\}
=P\subseteq\left\{  1,2,\ldots,n\right\}  $ (since $P$ is a subset of
$\left\{  1,2,\ldots,n\right\}  $). In other words, $i_{1},i_{2},\ldots,i_{k}$
are elements of $\left\{  1,2,\ldots,n\right\}  $. The same argument (but
applied to $m$, $Q$ and $\left(  j_{1},j_{2},\ldots,j_{k}\right)  $ instead of
$n$, $P$ and $\left(  i_{1},i_{2},\ldots,i_{k}\right)  $) yields that
$j_{1},j_{2},\ldots,j_{k}$ are elements of $\left\{  1,2,\ldots,m\right\}  $.
Hence, Corollary \ref{cor.adj(AB).cauchy-binet-general} (applied to $u=k$)
yields%
\begin{align}
&  \det\left(  \operatorname*{sub}\nolimits_{i_{1},i_{2},\ldots,i_{k}}%
^{j_{1},j_{2},\ldots,j_{k}}\left(  AB\right)  \right) \nonumber\\
&  =\sum_{1\leq g_{1}<g_{2}<\cdots<g_{k}\leq p}\det\left(  \operatorname*{sub}%
\nolimits_{i_{1},i_{2},\ldots,i_{k}}^{g_{1},g_{2},\ldots,g_{k}}A\right)
\cdot\det\left(  \operatorname*{sub}\nolimits_{g_{1},g_{2},\ldots,g_{k}%
}^{j_{1},j_{2},\ldots,j_{k}}B\right)  .
\label{pf.cor.sol.addexe.jacobi-complement.CB.1}%
\end{align}
(Here, the summation sign \textquotedblleft$\sum_{1\leq g_{1}<g_{2}%
<\cdots<g_{k}\leq p}$\textquotedblright\ has to be interpreted as
\textquotedblleft$\sum_{\substack{\left(  g_{1},g_{2},\ldots,g_{k}\right)
\in\left\{  1,2,\ldots,p\right\}  ^{k};\\g_{1}<g_{2}<\cdots<g_{k}}%
}$\textquotedblright, in analogy to Remark \ref{rmk.cauchy-binet.sumsign}.)

Now, let $S$ denote the set $\left\{  1,2,\ldots,p\right\}  $. Clearly, $S$ is
a set of integers. Recall that $\mathcal{P}_{k}\left(  S\right)  $ denotes the
set of all $k$-element subsets of $S$. In other words,
\[
\mathcal{P}_{k}\left(  S\right)  =\left\{  T\subseteq S\ \mid\ \left\vert
T\right\vert =k\right\}  .
\]

Let $\mathbf{I}$ denote the set
\[
\left\{  \left(  g_{1},g_{2},\ldots,g_{k}\right)  \in S^{k}\ \mid\ g_{1}%
<g_{2}<\cdots<g_{k}\right\}  .
\]
Thus, every element of $\mathbf{I}$ has the form $\left(  g_{1},g_{2}%
,\ldots,g_{k}\right)  $ for some $\left(  g_{1},g_{2},\ldots,g_{k}\right)  \in
S^{k}$.

Lemma \ref{lem.sol.addexe.jacobi-complement.increase-bij} (applied to $u=k$)
shows that the map%
\begin{align*}
\mathcal{P}_{k}\left(  S\right)   &  \rightarrow\mathbf{I},\\
R  &  \mapsto w\left(  R\right)
\end{align*}
is well-defined and a bijection.

The definition of $\mathbf{I}$ yields
\[
\mathbf{I}=\left\{  \left(  g_{1},g_{2},\ldots,g_{k}\right)  \in S^{k}%
\ \mid\ g_{1}<g_{2}<\cdots<g_{k}\right\}  .
\]
Hence, we have the following equality of summation signs:%
\[
\sum_{\left(  g_{1},g_{2},\ldots,g_{k}\right)  \in\mathbf{I}}=\sum
_{\substack{\left(  g_{1},g_{2},\ldots,g_{k}\right)  \in S^{k};\\g_{1}%
<g_{2}<\cdots<g_{k}}}=\sum_{\substack{\left(  g_{1},g_{2},\ldots,g_{k}\right)
\in\left\{  1,2,\ldots,p\right\}  ^{k};\\g_{1}<g_{2}<\cdots<g_{k}}}
\]
(since $S=\left\{  1,2,\ldots,p\right\}  $). Comparing this with%
\[
\sum_{1\leq g_{1}<g_{2}<\cdots<g_{k}\leq p}=\sum_{\substack{\left(
g_{1},g_{2},\ldots,g_{k}\right)  \in\left\{  1,2,\ldots,p\right\}
^{k};\\g_{1}<g_{2}<\cdots<g_{k}}},
\]
we obtain%
\[
\sum_{1\leq g_{1}<g_{2}<\cdots<g_{k}\leq p}=\sum_{\left(  g_{1},g_{2}%
,\ldots,g_{k}\right)  \in\mathbf{I}}.
\]
Now, (\ref{pf.cor.sol.addexe.jacobi-complement.CB.1}) becomes%
\begin{align*}
&  \det\left(  \operatorname*{sub}\nolimits_{i_{1},i_{2},\ldots,i_{k}}%
^{j_{1},j_{2},\ldots,j_{k}}\left(  AB\right)  \right) \\
&  =\underbrace{\sum_{1\leq g_{1}<g_{2}<\cdots<g_{k}\leq p}}_{=\sum_{\left(
g_{1},g_{2},\ldots,g_{k}\right)  \in\mathbf{I}}}\det\left(
\underbrace{\operatorname*{sub}\nolimits_{i_{1},i_{2},\ldots,i_{k}}%
^{g_{1},g_{2},\ldots,g_{k}}A}_{=\operatorname*{sub}\nolimits_{\left(
i_{1},i_{2},\ldots,i_{k}\right)  }^{\left(  g_{1},g_{2},\ldots,g_{k}\right)
}A}\right)  \cdot\det\left(  \underbrace{\operatorname*{sub}\nolimits_{g_{1}%
,g_{2},\ldots,g_{k}}^{j_{1},j_{2},\ldots,j_{k}}B}_{=\operatorname*{sub}%
\nolimits_{\left(  g_{1},g_{2},\ldots,g_{k}\right)  }^{\left(  j_{1}%
,j_{2},\ldots,j_{k}\right)  }B}\right) \\
&  =\sum_{\left(  g_{1},g_{2},\ldots,g_{k}\right)  \in\mathbf{I}}\det\left(
\operatorname*{sub}\nolimits_{\left(  i_{1},i_{2},\ldots,i_{k}\right)
}^{\left(  g_{1},g_{2},\ldots,g_{k}\right)  }A\right)  \cdot\det\left(
\operatorname*{sub}\nolimits_{\left(  g_{1},g_{2},\ldots,g_{k}\right)
}^{\left(  j_{1},j_{2},\ldots,j_{k}\right)  }B\right) \\
&  =\underbrace{\sum_{R\in\mathcal{P}_{k}\left(  S\right)  }}_{\substack{=\sum
_{R\in\left\{  T\subseteq S\ \mid\ \left\vert T\right\vert =k\right\}
}\\\text{(since }\mathcal{P}_{k}\left(  S\right)  =\left\{  T\subseteq
S\ \mid\ \left\vert T\right\vert =k\right\}  \text{)}}}\det\left(
\underbrace{\operatorname*{sub}\nolimits_{\left(  i_{1},i_{2},\ldots
,i_{k}\right)  }^{w\left(  R\right)  }A}_{\substack{=\operatorname*{sub}%
\nolimits_{w\left(  P\right)  }^{w\left(  R\right)  }A\\\text{(since }\left(
i_{1},i_{2},\ldots,i_{k}\right)  =w\left(  P\right)  \text{)}}}\right)
\cdot\det\left(  \underbrace{\operatorname*{sub}\nolimits_{w\left(  R\right)
}^{\left(  j_{1},j_{2},\ldots,j_{k}\right)  }B}%
_{\substack{=\operatorname*{sub}\nolimits_{w\left(  R\right)  }^{w\left(
Q\right)  }B\\\text{(since }\left(  j_{1},j_{2},\ldots,j_{k}\right)  =w\left(
Q\right)  \text{)}}}\right) \\
&  \ \ \ \ \ \ \ \ \ \ \left(
\begin{array}
[c]{c}%
\text{here, we have substituted }w\left(  R\right)  \text{ for }\left(
g_{1},g_{2},\ldots,g_{k}\right)  \text{, since}\\
\text{the map }\mathcal{P}_{k}\left(  S\right)  \rightarrow\mathbf{I}%
,\ R\mapsto w\left(  R\right)  \text{ is a bijection}%
\end{array}
\right) \\
&  =\underbrace{\sum_{R\in\left\{  T\subseteq S\ \mid\ \left\vert T\right\vert
=k\right\}  }}_{\substack{=\sum_{\substack{R\subseteq S;\\\left\vert
R\right\vert =k}}=\sum_{\substack{R\subseteq\left\{  1,2,\ldots,p\right\}
;\\\left\vert R\right\vert =k}}\\\text{(since }S=\left\{  1,2,\ldots
,p\right\}  \text{)}}}\det\left(  \operatorname*{sub}\nolimits_{w\left(
P\right)  }^{w\left(  R\right)  }A\right)  \cdot\det\left(
\operatorname*{sub}\nolimits_{w\left(  R\right)  }^{w\left(  Q\right)
}B\right) \\
&  =\sum_{\substack{R\subseteq\left\{  1,2,\ldots,p\right\}  ;\\\left\vert
R\right\vert =k}}\det\left(  \operatorname*{sub}\nolimits_{w\left(  P\right)
}^{w\left(  R\right)  }A\right)  \cdot\det\left(  \operatorname*{sub}%
\nolimits_{w\left(  R\right)  }^{w\left(  Q\right)  }B\right)  .
\end{align*}
Comparing this with%
\begin{align*}
\det\left(  \underbrace{\operatorname*{sub}\nolimits_{i_{1},i_{2},\ldots
,i_{k}}^{j_{1},j_{2},\ldots,j_{k}}\left(  AB\right)  }_{=\operatorname*{sub}%
\nolimits_{\left(  i_{1},i_{2},\ldots,i_{k}\right)  }^{\left(  j_{1}%
,j_{2},\ldots,j_{k}\right)  }\left(  AB\right)  }\right)   &  =\det\left(
\operatorname*{sub}\nolimits_{\left(  i_{1},i_{2},\ldots,i_{k}\right)
}^{\left(  j_{1},j_{2},\ldots,j_{k}\right)  }\left(  AB\right)  \right)
=\det\left(  \operatorname*{sub}\nolimits_{w\left(  P\right)  }^{w\left(
Q\right)  }\left(  AB\right)  \right) \\
&  \ \ \ \ \ \ \ \ \ \ \left(
\begin{array}
[c]{c}%
\text{since }\left(  i_{1},i_{2},\ldots,i_{k}\right)  =w\left(  P\right) \\
\text{and }\left(  j_{1},j_{2},\ldots,j_{k}\right)  =w\left(  Q\right)
\end{array}
\right)  ,
\end{align*}
we obtain%
\[
\det\left(  \operatorname*{sub}\nolimits_{w\left(  P\right)  }^{w\left(
Q\right)  }\left(  AB\right)  \right)  =\sum_{\substack{R\subseteq\left\{
1,2,\ldots,p\right\}  ;\\\left\vert R\right\vert =k}}\det\left(
\operatorname*{sub}\nolimits_{w\left(  P\right)  }^{w\left(  R\right)
}A\right)  \cdot\det\left(  \operatorname*{sub}\nolimits_{w\left(  R\right)
}^{w\left(  Q\right)  }B\right)  .
\]
This proves Corollary \ref{cor.sol.addexe.jacobi-complement.CB}.
\end{proof}
\end{verlong}

Our next step towards solving Exercise \ref{addexe.jacobi-complement} again is
the following fact:

\begin{proposition}
\label{prop.sol.addexe.jacobi-complement.genform}Let $n\in\mathbb{N}$ and
$m\in\mathbb{N}$. For any subset $I$ of $\left\{  1,2,\ldots,n\right\}  $, we
let $\widetilde{I}$ denote the complement $\left\{  1,2,\ldots,n\right\}
\setminus I$ of $I$.

Let $A$ be an $n\times n$-matrix. Let $B$ be an $n\times m$-matrix. Let $P$ be
a subset of $\left\{  1,2,\ldots,n\right\}  $. Let $Q$ be a subset of
$\left\{  1,2,\ldots,m\right\}  $ such that $\left\vert P\right\vert
=\left\vert Q\right\vert $. Then,%
\begin{align*}
&  \det A\cdot\det\left(  \operatorname*{sub}\nolimits_{w\left(  P\right)
}^{w\left(  Q\right)  }B\right) \\
&  =\sum_{\substack{K\subseteq\left\{  1,2,\ldots,n\right\}  ;\\\left\vert
K\right\vert =\left\vert P\right\vert }}\left(  -1\right)  ^{\sum P+\sum
K}\det\left(  \operatorname*{sub}\nolimits_{w\left(  \widetilde{K}\right)
}^{w\left(  \widetilde{P}\right)  }A\right)  \det\left(  \operatorname*{sub}%
\nolimits_{w\left(  K\right)  }^{w\left(  Q\right)  }\left(  AB\right)
\right)  .
\end{align*}

\end{proposition}

\begin{vershort}
\begin{proof}
[Proof of Proposition \ref{prop.sol.addexe.jacobi-complement.genform}.]Set
$k=\left\vert P\right\vert $. Clearly, $k=\left\vert P\right\vert =\left\vert
Q\right\vert $.

If $K$ is a subset of $\left\{  1,2,\ldots,n\right\}  $ satisfying $\left\vert
K\right\vert =\left\vert P\right\vert $, then%
\begin{align}
&  \det\left(  \operatorname*{sub}\nolimits_{w\left(  K\right)  }^{w\left(
Q\right)  }\left(  AB\right)  \right) \nonumber\\
&  =\sum_{\substack{R\subseteq\left\{  1,2,\ldots,n\right\}  ;\\\left\vert
R\right\vert =k}}\det\left(  \operatorname*{sub}\nolimits_{w\left(  K\right)
}^{w\left(  R\right)  }A\right)  \cdot\det\left(  \operatorname*{sub}%
\nolimits_{w\left(  R\right)  }^{w\left(  Q\right)  }B\right)
\label{pf.prop.sol.addexe.jacobi-complement.genform.short.1}%
\end{align}
\footnote{\textit{Proof of
(\ref{pf.prop.sol.addexe.jacobi-complement.genform.short.1}):} Let $K$ be a
subset of $\left\{  1,2,\ldots,n\right\}  $ satisfying $\left\vert
K\right\vert =\left\vert P\right\vert $. Then, $\left\vert K\right\vert
=\left\vert P\right\vert =k$ (since $k=\left\vert P\right\vert $). Also,
$\left\vert Q\right\vert =k$ (since $k=\left\vert Q\right\vert $). Hence,
Corollary \ref{cor.sol.addexe.jacobi-complement.CB} (applied to $n$ and $K$
instead of $p$ and $P$) yields
\[
\det\left(  \operatorname*{sub}\nolimits_{w\left(  K\right)  }^{w\left(
Q\right)  }\left(  AB\right)  \right)  =\sum_{\substack{R\subseteq\left\{
1,2,\ldots,n\right\}  ;\\\left\vert R\right\vert =k}}\det\left(
\operatorname*{sub}\nolimits_{w\left(  K\right)  }^{w\left(  R\right)
}A\right)  \cdot\det\left(  \operatorname*{sub}\nolimits_{w\left(  R\right)
}^{w\left(  Q\right)  }B\right)  .
\]
This proves (\ref{pf.prop.sol.addexe.jacobi-complement.genform.short.1}).}.

Now,%
\begin{align}
&  \sum_{\substack{K\subseteq\left\{  1,2,\ldots,n\right\}  ;\\\left\vert
K\right\vert =\left\vert P\right\vert }}\left(  -1\right)  ^{\sum P+\sum
K}\det\left(  \operatorname*{sub}\nolimits_{w\left(  \widetilde{K}\right)
}^{w\left(  \widetilde{P}\right)  }A\right)  \underbrace{\det\left(
\operatorname*{sub}\nolimits_{w\left(  K\right)  }^{w\left(  Q\right)
}\left(  AB\right)  \right)  }_{\substack{=\sum_{\substack{R\subseteq\left\{
1,2,\ldots,n\right\}  ;\\\left\vert R\right\vert =k}}\det\left(
\operatorname*{sub}\nolimits_{w\left(  K\right)  }^{w\left(  R\right)
}A\right)  \cdot\det\left(  \operatorname*{sub}\nolimits_{w\left(  R\right)
}^{w\left(  Q\right)  }B\right)  \\\text{(by
(\ref{pf.prop.sol.addexe.jacobi-complement.genform.short.1}))}}}\nonumber\\
&  =\sum_{\substack{K\subseteq\left\{  1,2,\ldots,n\right\}  ;\\\left\vert
K\right\vert =\left\vert P\right\vert }}\left(  -1\right)  ^{\sum P+\sum
K}\det\left(  \operatorname*{sub}\nolimits_{w\left(  \widetilde{K}\right)
}^{w\left(  \widetilde{P}\right)  }A\right) \nonumber\\
&  \ \ \ \ \ \ \ \ \ \ \left(  \sum_{\substack{R\subseteq\left\{
1,2,\ldots,n\right\}  ;\\\left\vert R\right\vert =k}}\det\left(
\operatorname*{sub}\nolimits_{w\left(  K\right)  }^{w\left(  R\right)
}A\right)  \cdot\det\left(  \operatorname*{sub}\nolimits_{w\left(  R\right)
}^{w\left(  Q\right)  }B\right)  \right) \nonumber\\
&  =\underbrace{\sum_{\substack{K\subseteq\left\{  1,2,\ldots,n\right\}
;\\\left\vert K\right\vert =\left\vert P\right\vert }}\ \ \sum
_{\substack{R\subseteq\left\{  1,2,\ldots,n\right\}  ;\\\left\vert
R\right\vert =k}}}_{=\sum_{\substack{R\subseteq\left\{  1,2,\ldots,n\right\}
;\\\left\vert R\right\vert =k}}\ \ \sum_{\substack{K\subseteq\left\{
1,2,\ldots,n\right\}  ;\\\left\vert K\right\vert =\left\vert P\right\vert }%
}}\underbrace{\left(  -1\right)  ^{\sum P+\sum K}}_{\substack{=\left(
-1\right)  ^{\sum K+\sum P}\\\text{(since }\sum P+\sum K=\sum K+\sum
P\text{)}}}\nonumber\\
&  \ \ \ \ \ \ \ \ \ \ \underbrace{\det\left(  \operatorname*{sub}%
\nolimits_{w\left(  \widetilde{K}\right)  }^{w\left(  \widetilde{P}\right)
}A\right)  \det\left(  \operatorname*{sub}\nolimits_{w\left(  K\right)
}^{w\left(  R\right)  }A\right)  }_{=\det\left(  \operatorname*{sub}%
\nolimits_{w\left(  K\right)  }^{w\left(  R\right)  }A\right)  \det\left(
\operatorname*{sub}\nolimits_{w\left(  \widetilde{K}\right)  }^{w\left(
\widetilde{P}\right)  }A\right)  }\cdot\det\left(  \operatorname*{sub}%
\nolimits_{w\left(  R\right)  }^{w\left(  Q\right)  }B\right) \nonumber\\
&  =\sum_{\substack{R\subseteq\left\{  1,2,\ldots,n\right\}  ;\\\left\vert
R\right\vert =k}}\ \ \sum_{\substack{K\subseteq\left\{  1,2,\ldots,n\right\}
;\\\left\vert K\right\vert =\left\vert P\right\vert }}\left(  -1\right)
^{\sum K+\sum P}\nonumber\\
&  \ \ \ \ \ \ \ \ \ \ \det\left(  \operatorname*{sub}\nolimits_{w\left(
K\right)  }^{w\left(  R\right)  }A\right)  \det\left(  \operatorname*{sub}%
\nolimits_{w\left(  \widetilde{K}\right)  }^{w\left(  \widetilde{P}\right)
}A\right)  \cdot\det\left(  \operatorname*{sub}\nolimits_{w\left(  R\right)
}^{w\left(  Q\right)  }B\right) \nonumber\\
&  =\sum_{\substack{R\subseteq\left\{  1,2,\ldots,n\right\}  ;\\\left\vert
R\right\vert =k}}\left(  \sum_{\substack{K\subseteq\left\{  1,2,\ldots
,n\right\}  ;\\\left\vert K\right\vert =\left\vert P\right\vert }}\left(
-1\right)  ^{\sum K+\sum P}\det\left(  \operatorname*{sub}\nolimits_{w\left(
K\right)  }^{w\left(  R\right)  }A\right)  \det\left(  \operatorname*{sub}%
\nolimits_{w\left(  \widetilde{K}\right)  }^{w\left(  \widetilde{P}\right)
}A\right)  \right) \nonumber\\
&  \ \ \ \ \ \ \ \ \ \ \cdot\det\left(  \operatorname*{sub}\nolimits_{w\left(
R\right)  }^{w\left(  Q\right)  }B\right)  .
\label{pf.prop.sol.addexe.jacobi-complement.genform.short.3}%
\end{align}

On the other hand, for every subset $G$ of $\left\{  1,2,\ldots,n\right\}  $,
we have%
\begin{equation}
\det A=\sum_{\substack{K\subseteq\left\{  1,2,\ldots,n\right\}  ;\\\left\vert
K\right\vert =\left\vert G\right\vert }}\left(  -1\right)  ^{\sum K+\sum
G}\det\left(  \operatorname*{sub}\nolimits_{w\left(  K\right)  }^{w\left(
G\right)  }A\right)  \det\left(  \operatorname*{sub}\nolimits_{w\left(
\widetilde{K}\right)  }^{w\left(  \widetilde{G}\right)  }A\right)  .
\label{pf.prop.sol.addexe.jacobi-complement.genform.short.6}%
\end{equation}
(Indeed, this is precisely the claim of Theorem \ref{thm.det.laplace-multi}
\textbf{(b)}, with the variables $P$ and $Q$ renamed as $K$ and $G$.) Applying
(\ref{pf.prop.sol.addexe.jacobi-complement.genform.short.6}) to $G=P$, we
obtain%
\begin{equation}
\det A=\sum_{\substack{K\subseteq\left\{  1,2,\ldots,n\right\}  ;\\\left\vert
K\right\vert =\left\vert P\right\vert }}\left(  -1\right)  ^{\sum K+\sum
P}\det\left(  \operatorname*{sub}\nolimits_{w\left(  K\right)  }^{w\left(
P\right)  }A\right)  \det\left(  \operatorname*{sub}\nolimits_{w\left(
\widetilde{K}\right)  }^{w\left(  \widetilde{P}\right)  }A\right)  .
\label{pf.prop.sol.addexe.jacobi-complement.genform.short.6b}%
\end{equation}

For any two objects $i$ and $j$, we define $\delta_{i,j}$ to be the element $%
\begin{cases}
1, & \text{if }i=j;\\
0, & \text{if }i\neq j
\end{cases}
$ of $\mathbb{K}$.

If $R$ is a subset of $\left\{  1,2,\ldots,n\right\}  $ satisfying $\left\vert
R\right\vert =k$, then%
\begin{align}
&  \sum_{\substack{K\subseteq\left\{  1,2,\ldots,n\right\}  ;\\\left\vert
K\right\vert =\left\vert P\right\vert }}\left(  -1\right)  ^{\sum K+\sum
P}\det\left(  \operatorname*{sub}\nolimits_{w\left(  K\right)  }^{w\left(
R\right)  }A\right)  \det\left(  \operatorname*{sub}\nolimits_{w\left(
\widetilde{K}\right)  }^{w\left(  \widetilde{P}\right)  }A\right) \nonumber\\
&  =\delta_{R,P}\det A
\label{pf.prop.sol.addexe.jacobi-complement.genform.short.5}%
\end{align}
\footnote{\textit{Proof of
(\ref{pf.prop.sol.addexe.jacobi-complement.genform.short.5}):} Let $R$ be a
subset of $\left\{  1,2,\ldots,n\right\}  $ satisfying $\left\vert
R\right\vert =k$. We must prove
(\ref{pf.prop.sol.addexe.jacobi-complement.genform.short.5}).
\par
We have
\[
\underbrace{\delta_{P,P}}_{\substack{=1\\\text{(since }P=P\text{)}}}\det
A=\det A=\sum_{\substack{K\subseteq\left\{  1,2,\ldots,n\right\}
;\\\left\vert K\right\vert =\left\vert P\right\vert }}\left(  -1\right)
^{\sum K+\sum P}\det\left(  \operatorname*{sub}\nolimits_{w\left(  K\right)
}^{w\left(  P\right)  }A\right)  \det\left(  \operatorname*{sub}%
\nolimits_{w\left(  \widetilde{K}\right)  }^{w\left(  \widetilde{P}\right)
}A\right)
\]
(by (\ref{pf.prop.sol.addexe.jacobi-complement.genform.short.6}), applied to
$G=P$). In other words,
(\ref{pf.prop.sol.addexe.jacobi-complement.genform.short.5}) holds if $R=P$.
Hence, for the rest of our proof of
(\ref{pf.prop.sol.addexe.jacobi-complement.genform.short.5}), we WLOG assume
that $R\neq P$. Thus, $\delta_{R,P}=0$.
\par
We have $R\neq P$ and thus $P\neq R$. Also, $\left\vert P\right\vert
=k=\left\vert R\right\vert $.
\par
For every subset $G$ of $\left\{  1,2,\ldots,n\right\}  $ satisfying
$\left\vert G\right\vert =\left\vert R\right\vert $ and $G\neq R$, we have%
\begin{equation}
0=\sum_{\substack{K\subseteq\left\{  1,2,\ldots,n\right\}  ;\\\left\vert
K\right\vert =\left\vert G\right\vert }}\left(  -1\right)  ^{\sum K+\sum
G}\det\left(  \operatorname*{sub}\nolimits_{w\left(  K\right)  }^{w\left(
R\right)  }A\right)  \det\left(  \operatorname*{sub}\nolimits_{w\left(
\widetilde{K}\right)  }^{w\left(  \widetilde{G}\right)  }A\right)  .
\label{pf.prop.sol.addexe.jacobi-complement.genform.short.7}%
\end{equation}
(Indeed, this is precisely the claim of Exercise \ref{exe.det.laplace-multi.0}
\textbf{(b)}, with the variables $Q$ and $P$ renamed as $G$ and $K$.) Applying
(\ref{pf.prop.sol.addexe.jacobi-complement.genform.short.7}) to $G=P$, we
obtain%
\[
0=\sum_{\substack{K\subseteq\left\{  1,2,\ldots,n\right\}  ;\\\left\vert
K\right\vert =\left\vert P\right\vert }}\left(  -1\right)  ^{\sum K+\sum
P}\det\left(  \operatorname*{sub}\nolimits_{w\left(  K\right)  }^{w\left(
R\right)  }A\right)  \det\left(  \operatorname*{sub}\nolimits_{w\left(
\widetilde{K}\right)  }^{w\left(  \widetilde{P}\right)  }A\right)
\]
(since $\left\vert P\right\vert =\left\vert R\right\vert $ and $P\neq R$).
Hence,%
\begin{align*}
&  \sum_{\substack{K\subseteq\left\{  1,2,\ldots,n\right\}  ;\\\left\vert
K\right\vert =\left\vert P\right\vert }}\left(  -1\right)  ^{\sum K+\sum
P}\det\left(  \operatorname*{sub}\nolimits_{w\left(  K\right)  }^{w\left(
R\right)  }A\right)  \det\left(  \operatorname*{sub}\nolimits_{w\left(
\widetilde{K}\right)  }^{w\left(  \widetilde{P}\right)  }A\right) \\
&  =0=\underbrace{0}_{=\delta_{R,P}}\det A=\delta_{R,P}\det A.
\end{align*}
This proves (\ref{pf.prop.sol.addexe.jacobi-complement.genform.short.5}).}.

Now, (\ref{pf.prop.sol.addexe.jacobi-complement.genform.short.3}) becomes%
\begin{align}
&  \sum_{\substack{K\subseteq\left\{  1,2,\ldots,n\right\}  ;\\\left\vert
K\right\vert =\left\vert P\right\vert }}\left(  -1\right)  ^{\sum P+\sum
K}\det\left(  \operatorname*{sub}\nolimits_{w\left(  \widetilde{K}\right)
}^{w\left(  \widetilde{P}\right)  }A\right)  \det\left(  \operatorname*{sub}%
\nolimits_{w\left(  K\right)  }^{w\left(  Q\right)  }\left(  AB\right)
\right) \nonumber\\
&  =\sum_{\substack{R\subseteq\left\{  1,2,\ldots,n\right\}  ;\\\left\vert
R\right\vert =k}}\underbrace{\left(  \sum_{\substack{K\subseteq\left\{
1,2,\ldots,n\right\}  ;\\\left\vert K\right\vert =\left\vert P\right\vert
}}\left(  -1\right)  ^{\sum K+\sum P}\det\left(  \operatorname*{sub}%
\nolimits_{w\left(  K\right)  }^{w\left(  R\right)  }A\right)  \det\left(
\operatorname*{sub}\nolimits_{w\left(  \widetilde{K}\right)  }^{w\left(
\widetilde{P}\right)  }A\right)  \right)  }_{\substack{=\delta_{R,P}\det
A\\\text{(by (\ref{pf.prop.sol.addexe.jacobi-complement.genform.short.5}))}%
}}\nonumber\\
&  \ \ \ \ \ \ \ \ \ \ \cdot\det\left(  \operatorname*{sub}\nolimits_{w\left(
R\right)  }^{w\left(  Q\right)  }B\right) \nonumber\\
&  =\sum_{\substack{R\subseteq\left\{  1,2,\ldots,n\right\}  ;\\\left\vert
R\right\vert =k}}\delta_{R,P}\det A\cdot\det\left(  \operatorname*{sub}%
\nolimits_{w\left(  R\right)  }^{w\left(  Q\right)  }B\right)  .
\label{pf.prop.sol.addexe.jacobi-complement.genform.short.11}%
\end{align}

But $P$ is a subset of $\left\{  1,2,\ldots,n\right\}  $ and satisfies
$\left\vert P\right\vert =k$. In other words, $P$ is a subset $R$ of $\left\{
1,2,\ldots,n\right\}  $ satisfying $\left\vert R\right\vert =k$. Thus, the sum
$\sum_{\substack{R\subseteq\left\{  1,2,\ldots,n\right\}  ;\\\left\vert
R\right\vert =k}}\delta_{R,P}\det A\cdot\det\left(  \operatorname*{sub}%
\nolimits_{w\left(  R\right)  }^{w\left(  Q\right)  }B\right)  $ has an addend
for $R=P$. If we split off this addend from this sum, then we obtain%
\begin{align*}
&  \sum_{\substack{R\subseteq\left\{  1,2,\ldots,n\right\}  ;\\\left\vert
R\right\vert =k}}\delta_{R,P}\det A\cdot\det\left(  \operatorname*{sub}%
\nolimits_{w\left(  R\right)  }^{w\left(  Q\right)  }B\right) \\
&  =\underbrace{\delta_{P,P}}_{\substack{=1\\\text{(since }P=P\text{)}}}\det
A\cdot\det\left(  \operatorname*{sub}\nolimits_{w\left(  P\right)  }^{w\left(
Q\right)  }B\right)  +\sum_{\substack{R\subseteq\left\{  1,2,\ldots,n\right\}
;\\\left\vert R\right\vert =k;\\R\neq P}}\underbrace{\delta_{R,P}%
}_{\substack{=0\\\text{(since }R\neq P\text{)}}}\det A\cdot\det\left(
\operatorname*{sub}\nolimits_{w\left(  R\right)  }^{w\left(  Q\right)
}B\right) \\
&  =\det A\cdot\det\left(  \operatorname*{sub}\nolimits_{w\left(  P\right)
}^{w\left(  Q\right)  }B\right)  +\underbrace{\sum_{\substack{R\subseteq
\left\{  1,2,\ldots,n\right\}  ;\\\left\vert R\right\vert =k;\\R\neq P}}0\det
A\cdot\det\left(  \operatorname*{sub}\nolimits_{w\left(  R\right)  }^{w\left(
Q\right)  }B\right)  }_{=0}\\
&  =\det A\cdot\det\left(  \operatorname*{sub}\nolimits_{w\left(  P\right)
}^{w\left(  Q\right)  }B\right)  .
\end{align*}
Hence, (\ref{pf.prop.sol.addexe.jacobi-complement.genform.short.11}) becomes%
\begin{align*}
&  \sum_{\substack{K\subseteq\left\{  1,2,\ldots,n\right\}  ;\\\left\vert
K\right\vert =\left\vert P\right\vert }}\left(  -1\right)  ^{\sum P+\sum
K}\det\left(  \operatorname*{sub}\nolimits_{w\left(  \widetilde{K}\right)
}^{w\left(  \widetilde{P}\right)  }A\right)  \det\left(  \operatorname*{sub}%
\nolimits_{w\left(  K\right)  }^{w\left(  Q\right)  }\left(  AB\right)
\right) \\
&  =\sum_{\substack{R\subseteq\left\{  1,2,\ldots,n\right\}  ;\\\left\vert
R\right\vert =k}}\delta_{R,P}\det A\cdot\det\left(  \operatorname*{sub}%
\nolimits_{w\left(  R\right)  }^{w\left(  Q\right)  }B\right)  =\det
A\cdot\det\left(  \operatorname*{sub}\nolimits_{w\left(  P\right)  }^{w\left(
Q\right)  }B\right)  .
\end{align*}
This proves Proposition \ref{prop.sol.addexe.jacobi-complement.genform}.
\end{proof}
\end{vershort}

\begin{verlong}
\begin{proof}
[Proof of Proposition \ref{prop.sol.addexe.jacobi-complement.genform}.]The set
$P$ is a subset of $\left\{  1,2,\ldots,n\right\}  $, and thus is finite
(since $\left\{  1,2,\ldots,n\right\}  $ is finite). Thus, $\left\vert
P\right\vert \in\mathbb{N}$. Hence, we can define $k\in\mathbb{N}$ by
$k=\left\vert P\right\vert $. Consider this $k$. Clearly, $k=\left\vert
P\right\vert =\left\vert Q\right\vert $.

If $K$ is a subset of $\left\{  1,2,\ldots,n\right\}  $ satisfying $\left\vert
K\right\vert =\left\vert P\right\vert $, then%
\begin{align}
&  \det\left(  \operatorname*{sub}\nolimits_{w\left(  K\right)  }^{w\left(
Q\right)  }\left(  AB\right)  \right) \nonumber\\
&  =\sum_{\substack{R\subseteq\left\{  1,2,\ldots,n\right\}  ;\\\left\vert
R\right\vert =k}}\det\left(  \operatorname*{sub}\nolimits_{w\left(  K\right)
}^{w\left(  R\right)  }A\right)  \cdot\det\left(  \operatorname*{sub}%
\nolimits_{w\left(  R\right)  }^{w\left(  Q\right)  }B\right)
\label{pf.prop.sol.addexe.jacobi-complement.genform.1}%
\end{align}
\footnote{\textit{Proof of
(\ref{pf.prop.sol.addexe.jacobi-complement.genform.1}):} Let $K$ be a subset
of $\left\{  1,2,\ldots,n\right\}  $ satisfying $\left\vert K\right\vert
=\left\vert P\right\vert $. Then, $\left\vert K\right\vert =\left\vert
P\right\vert =k$ (since $k=\left\vert P\right\vert $). Also, $\left\vert
Q\right\vert =k$ (since $k=\left\vert Q\right\vert $). Hence, Corollary
\ref{cor.sol.addexe.jacobi-complement.CB} (applied to $n$ and $K$ instead of
$p$ and $P$) yields
\[
\det\left(  \operatorname*{sub}\nolimits_{w\left(  K\right)  }^{w\left(
Q\right)  }\left(  AB\right)  \right)  =\sum_{\substack{R\subseteq\left\{
1,2,\ldots,n\right\}  ;\\\left\vert R\right\vert =k}}\det\left(
\operatorname*{sub}\nolimits_{w\left(  K\right)  }^{w\left(  R\right)
}A\right)  \cdot\det\left(  \operatorname*{sub}\nolimits_{w\left(  R\right)
}^{w\left(  Q\right)  }B\right)  .
\]
This proves (\ref{pf.prop.sol.addexe.jacobi-complement.genform.1}).}.

Now,%
\begin{align}
&  \sum_{\substack{K\subseteq\left\{  1,2,\ldots,n\right\}  ;\\\left\vert
K\right\vert =\left\vert P\right\vert }}\left(  -1\right)  ^{\sum P+\sum
K}\det\left(  \operatorname*{sub}\nolimits_{w\left(  \widetilde{K}\right)
}^{w\left(  \widetilde{P}\right)  }A\right)  \underbrace{\det\left(
\operatorname*{sub}\nolimits_{w\left(  K\right)  }^{w\left(  Q\right)
}\left(  AB\right)  \right)  }_{\substack{=\sum_{\substack{R\subseteq\left\{
1,2,\ldots,n\right\}  ;\\\left\vert R\right\vert =k}}\det\left(
\operatorname*{sub}\nolimits_{w\left(  K\right)  }^{w\left(  R\right)
}A\right)  \cdot\det\left(  \operatorname*{sub}\nolimits_{w\left(  R\right)
}^{w\left(  Q\right)  }B\right)  \\\text{(by
(\ref{pf.prop.sol.addexe.jacobi-complement.genform.1}))}}}\nonumber\\
&  =\sum_{\substack{K\subseteq\left\{  1,2,\ldots,n\right\}  ;\\\left\vert
K\right\vert =\left\vert P\right\vert }}\left(  -1\right)  ^{\sum P+\sum
K}\det\left(  \operatorname*{sub}\nolimits_{w\left(  \widetilde{K}\right)
}^{w\left(  \widetilde{P}\right)  }A\right) \nonumber\\
&  \ \ \ \ \ \ \ \ \ \ \left(  \sum_{\substack{R\subseteq\left\{
1,2,\ldots,n\right\}  ;\\\left\vert R\right\vert =k}}\det\left(
\operatorname*{sub}\nolimits_{w\left(  K\right)  }^{w\left(  R\right)
}A\right)  \cdot\det\left(  \operatorname*{sub}\nolimits_{w\left(  R\right)
}^{w\left(  Q\right)  }B\right)  \right) \nonumber\\
&  =\underbrace{\sum_{\substack{K\subseteq\left\{  1,2,\ldots,n\right\}
;\\\left\vert K\right\vert =\left\vert P\right\vert }}\sum
_{\substack{R\subseteq\left\{  1,2,\ldots,n\right\}  ;\\\left\vert
R\right\vert =k}}}_{=\sum_{\substack{R\subseteq\left\{  1,2,\ldots,n\right\}
;\\\left\vert R\right\vert =k}}\sum_{\substack{K\subseteq\left\{
1,2,\ldots,n\right\}  ;\\\left\vert K\right\vert =\left\vert P\right\vert }%
}}\underbrace{\left(  -1\right)  ^{\sum P+\sum K}}_{\substack{=\left(
-1\right)  ^{\sum K+\sum P}\\\text{(since }\sum P+\sum K=\sum K+\sum
P\text{)}}}\nonumber\\
&  \ \ \ \ \ \ \ \ \ \ \underbrace{\det\left(  \operatorname*{sub}%
\nolimits_{w\left(  \widetilde{K}\right)  }^{w\left(  \widetilde{P}\right)
}A\right)  \det\left(  \operatorname*{sub}\nolimits_{w\left(  K\right)
}^{w\left(  R\right)  }A\right)  }_{=\det\left(  \operatorname*{sub}%
\nolimits_{w\left(  K\right)  }^{w\left(  R\right)  }A\right)  \det\left(
\operatorname*{sub}\nolimits_{w\left(  \widetilde{K}\right)  }^{w\left(
\widetilde{P}\right)  }A\right)  }\cdot\det\left(  \operatorname*{sub}%
\nolimits_{w\left(  R\right)  }^{w\left(  Q\right)  }B\right) \nonumber\\
&  =\sum_{\substack{R\subseteq\left\{  1,2,\ldots,n\right\}  ;\\\left\vert
R\right\vert =k}}\sum_{\substack{K\subseteq\left\{  1,2,\ldots,n\right\}
;\\\left\vert K\right\vert =\left\vert P\right\vert }}\left(  -1\right)
^{\sum K+\sum P}\nonumber\\
&  \ \ \ \ \ \ \ \ \ \ \det\left(  \operatorname*{sub}\nolimits_{w\left(
K\right)  }^{w\left(  R\right)  }A\right)  \det\left(  \operatorname*{sub}%
\nolimits_{w\left(  \widetilde{K}\right)  }^{w\left(  \widetilde{P}\right)
}A\right)  \cdot\det\left(  \operatorname*{sub}\nolimits_{w\left(  R\right)
}^{w\left(  Q\right)  }B\right) \nonumber\\
&  =\sum_{\substack{R\subseteq\left\{  1,2,\ldots,n\right\}  ;\\\left\vert
R\right\vert =k}}\left(  \sum_{\substack{K\subseteq\left\{  1,2,\ldots
,n\right\}  ;\\\left\vert K\right\vert =\left\vert P\right\vert }}\left(
-1\right)  ^{\sum K+\sum P}\det\left(  \operatorname*{sub}\nolimits_{w\left(
K\right)  }^{w\left(  R\right)  }A\right)  \det\left(  \operatorname*{sub}%
\nolimits_{w\left(  \widetilde{K}\right)  }^{w\left(  \widetilde{P}\right)
}A\right)  \right) \nonumber\\
&  \ \ \ \ \ \ \ \ \ \ \cdot\det\left(  \operatorname*{sub}\nolimits_{w\left(
R\right)  }^{w\left(  Q\right)  }B\right)  .
\label{pf.prop.sol.addexe.jacobi-complement.genform.3}%
\end{align}

On the other hand, for every subset $G$ of $\left\{  1,2,\ldots,n\right\}  $,
we have%
\begin{equation}
\det A=\sum_{\substack{K\subseteq\left\{  1,2,\ldots,n\right\}  ;\\\left\vert
K\right\vert =\left\vert G\right\vert }}\left(  -1\right)  ^{\sum K+\sum
G}\det\left(  \operatorname*{sub}\nolimits_{w\left(  K\right)  }^{w\left(
G\right)  }A\right)  \det\left(  \operatorname*{sub}\nolimits_{w\left(
\widetilde{K}\right)  }^{w\left(  \widetilde{G}\right)  }A\right)  .
\label{pf.prop.sol.addexe.jacobi-complement.genform.6}%
\end{equation}
(Indeed, this is precisely the claim of Theorem \ref{thm.det.laplace-multi}
\textbf{(b)}, with the variables $P$ and $Q$ renamed as $K$ and $G$.) Applying
(\ref{pf.prop.sol.addexe.jacobi-complement.genform.6}) to $G=P$, we obtain%
\begin{equation}
\det A=\sum_{\substack{K\subseteq\left\{  1,2,\ldots,n\right\}  ;\\\left\vert
K\right\vert =\left\vert P\right\vert }}\left(  -1\right)  ^{\sum K+\sum
P}\det\left(  \operatorname*{sub}\nolimits_{w\left(  K\right)  }^{w\left(
P\right)  }A\right)  \det\left(  \operatorname*{sub}\nolimits_{w\left(
\widetilde{K}\right)  }^{w\left(  \widetilde{P}\right)  }A\right)  .
\label{pf.prop.sol.addexe.jacobi-complement.genform.6b}%
\end{equation}

For any two objects $i$ and $j$, we define $\delta_{i,j}$ to be the element $%
\begin{cases}
1, & \text{if }i=j;\\
0, & \text{if }i\neq j
\end{cases}
$ of $\mathbb{K}$.

If $R$ is a subset of $\left\{  1,2,\ldots,n\right\}  $ satisfying $\left\vert
R\right\vert =k$, then%
\begin{align}
&  \sum_{\substack{K\subseteq\left\{  1,2,\ldots,n\right\}  ;\\\left\vert
K\right\vert =\left\vert P\right\vert }}\left(  -1\right)  ^{\sum K+\sum
P}\det\left(  \operatorname*{sub}\nolimits_{w\left(  K\right)  }^{w\left(
R\right)  }A\right)  \det\left(  \operatorname*{sub}\nolimits_{w\left(
\widetilde{K}\right)  }^{w\left(  \widetilde{P}\right)  }A\right) \nonumber\\
&  =\delta_{R,P}\det A \label{pf.prop.sol.addexe.jacobi-complement.genform.5}%
\end{align}
\footnote{\textit{Proof of
(\ref{pf.prop.sol.addexe.jacobi-complement.genform.5}):} Let $R$ be a subset
of $\left\{  1,2,\ldots,n\right\}  $ satisfying $\left\vert R\right\vert =k$.
We must prove (\ref{pf.prop.sol.addexe.jacobi-complement.genform.5}).
\par
We have
\begin{align*}
\underbrace{\delta_{P,P}}_{\substack{=1\\\text{(since }P=P\text{)}}}\det A  &
=\det A\\
&  =\sum_{\substack{K\subseteq\left\{  1,2,\ldots,n\right\}  ;\\\left\vert
K\right\vert =\left\vert P\right\vert }}\left(  -1\right)  ^{\sum K+\sum
P}\det\left(  \operatorname*{sub}\nolimits_{w\left(  K\right)  }^{w\left(
P\right)  }A\right)  \det\left(  \operatorname*{sub}\nolimits_{w\left(
\widetilde{K}\right)  }^{w\left(  \widetilde{P}\right)  }A\right)
\end{align*}
(by (\ref{pf.prop.sol.addexe.jacobi-complement.genform.6}), applied to $G=P$).
Thus,%
\[
\sum_{\substack{K\subseteq\left\{  1,2,\ldots,n\right\}  ;\\\left\vert
K\right\vert =\left\vert P\right\vert }}\left(  -1\right)  ^{\sum K+\sum
P}\det\left(  \operatorname*{sub}\nolimits_{w\left(  K\right)  }^{w\left(
P\right)  }A\right)  \det\left(  \operatorname*{sub}\nolimits_{w\left(
\widetilde{K}\right)  }^{w\left(  \widetilde{P}\right)  }A\right)
=\delta_{P,P}\det A.
\]
In other words, (\ref{pf.prop.sol.addexe.jacobi-complement.genform.5}) holds
if $R=P$. Hence, for the rest of our proof of
(\ref{pf.prop.sol.addexe.jacobi-complement.genform.5}), we can WLOG assume
that we don't have $R=P$. Assume this.
\par
We have $R\neq P$ (since we don't have $R=P$). Thus, $\delta_{R,P}=0$.
\par
We have $R\neq P$ and thus $P\neq R$. Also, $\left\vert P\right\vert
=k=\left\vert R\right\vert $ (since $\left\vert R\right\vert =k$).
\par
For every subset $G$ of $\left\{  1,2,\ldots,n\right\}  $ satisfying
$\left\vert G\right\vert =\left\vert R\right\vert $ and $G\neq R$, we have%
\begin{equation}
0=\sum_{\substack{K\subseteq\left\{  1,2,\ldots,n\right\}  ;\\\left\vert
K\right\vert =\left\vert G\right\vert }}\left(  -1\right)  ^{\sum K+\sum
G}\det\left(  \operatorname*{sub}\nolimits_{w\left(  K\right)  }^{w\left(
R\right)  }A\right)  \det\left(  \operatorname*{sub}\nolimits_{w\left(
\widetilde{K}\right)  }^{w\left(  \widetilde{G}\right)  }A\right)  .
\label{pf.prop.sol.addexe.jacobi-complement.genform.7}%
\end{equation}
(Indeed, this is precisely the claim of Exercise \ref{exe.det.laplace-multi.0}
\textbf{(b)}, with the variables $Q$ and $P$ renamed as $G$ and $K$.) Applying
(\ref{pf.prop.sol.addexe.jacobi-complement.genform.7}) to $G=P$, we obtain%
\[
0=\sum_{\substack{K\subseteq\left\{  1,2,\ldots,n\right\}  ;\\\left\vert
K\right\vert =\left\vert P\right\vert }}\left(  -1\right)  ^{\sum K+\sum
P}\det\left(  \operatorname*{sub}\nolimits_{w\left(  K\right)  }^{w\left(
R\right)  }A\right)  \det\left(  \operatorname*{sub}\nolimits_{w\left(
\widetilde{K}\right)  }^{w\left(  \widetilde{P}\right)  }A\right)
\]
(since $\left\vert P\right\vert =\left\vert R\right\vert $ and $P\neq R$).
Hence,%
\begin{align*}
&  \sum_{\substack{K\subseteq\left\{  1,2,\ldots,n\right\}  ;\\\left\vert
K\right\vert =\left\vert P\right\vert }}\left(  -1\right)  ^{\sum K+\sum
P}\det\left(  \operatorname*{sub}\nolimits_{w\left(  K\right)  }^{w\left(
R\right)  }A\right)  \det\left(  \operatorname*{sub}\nolimits_{w\left(
\widetilde{K}\right)  }^{w\left(  \widetilde{P}\right)  }A\right) \\
&  =0=\underbrace{0}_{=\delta_{R,P}}\det A=\delta_{R,P}\det A.
\end{align*}
This proves (\ref{pf.prop.sol.addexe.jacobi-complement.genform.5}).}.

Now, (\ref{pf.prop.sol.addexe.jacobi-complement.genform.3}) becomes%
\begin{align}
&  \sum_{\substack{K\subseteq\left\{  1,2,\ldots,n\right\}  ;\\\left\vert
K\right\vert =\left\vert P\right\vert }}\left(  -1\right)  ^{\sum P+\sum
K}\det\left(  \operatorname*{sub}\nolimits_{w\left(  \widetilde{K}\right)
}^{w\left(  \widetilde{P}\right)  }A\right)  \det\left(  \operatorname*{sub}%
\nolimits_{w\left(  K\right)  }^{w\left(  Q\right)  }\left(  AB\right)
\right) \nonumber\\
&  =\sum_{\substack{R\subseteq\left\{  1,2,\ldots,n\right\}  ;\\\left\vert
R\right\vert =k}}\underbrace{\left(  \sum_{\substack{K\subseteq\left\{
1,2,\ldots,n\right\}  ;\\\left\vert K\right\vert =\left\vert P\right\vert
}}\left(  -1\right)  ^{\sum K+\sum P}\det\left(  \operatorname*{sub}%
\nolimits_{w\left(  K\right)  }^{w\left(  R\right)  }A\right)  \det\left(
\operatorname*{sub}\nolimits_{w\left(  \widetilde{K}\right)  }^{w\left(
\widetilde{P}\right)  }A\right)  \right)  }_{\substack{=\delta_{R,P}\det
A\\\text{(by (\ref{pf.prop.sol.addexe.jacobi-complement.genform.5}))}%
}}\nonumber\\
&  \ \ \ \ \ \ \ \ \ \ \cdot\det\left(  \operatorname*{sub}\nolimits_{w\left(
R\right)  }^{w\left(  Q\right)  }B\right) \nonumber\\
&  =\sum_{\substack{R\subseteq\left\{  1,2,\ldots,n\right\}  ;\\\left\vert
R\right\vert =k}}\delta_{R,P}\det A\cdot\det\left(  \operatorname*{sub}%
\nolimits_{w\left(  R\right)  }^{w\left(  Q\right)  }B\right)  .
\label{pf.prop.sol.addexe.jacobi-complement.genform.11}%
\end{align}

But $P$ is a subset of $\left\{  1,2,\ldots,n\right\}  $ and satisfies
$\left\vert P\right\vert =k$. In other words, $P$ is a subset $R$ of $\left\{
1,2,\ldots,n\right\}  $ satisfying $\left\vert R\right\vert =k$. Thus, the sum
$\sum_{\substack{R\subseteq\left\{  1,2,\ldots,n\right\}  ;\\\left\vert
R\right\vert =k}}\delta_{R,P}\det A\cdot\det\left(  \operatorname*{sub}%
\nolimits_{w\left(  R\right)  }^{w\left(  Q\right)  }B\right)  $ has an addend
for $R=P$. If we split off this addend from this sum, then we obtain%
\begin{align*}
&  \sum_{\substack{R\subseteq\left\{  1,2,\ldots,n\right\}  ;\\\left\vert
R\right\vert =k}}\delta_{R,P}\det A\cdot\det\left(  \operatorname*{sub}%
\nolimits_{w\left(  R\right)  }^{w\left(  Q\right)  }B\right) \\
&  =\underbrace{\delta_{P,P}}_{\substack{=1\\\text{(since }P=P\text{)}}}\det
A\cdot\det\left(  \operatorname*{sub}\nolimits_{w\left(  P\right)  }^{w\left(
Q\right)  }B\right)  +\sum_{\substack{R\subseteq\left\{  1,2,\ldots,n\right\}
;\\\left\vert R\right\vert =k;\\R\neq P}}\underbrace{\delta_{R,P}%
}_{\substack{=0\\\text{(since }R\neq P\text{)}}}\det A\cdot\det\left(
\operatorname*{sub}\nolimits_{w\left(  R\right)  }^{w\left(  Q\right)
}B\right) \\
&  =\det A\cdot\det\left(  \operatorname*{sub}\nolimits_{w\left(  P\right)
}^{w\left(  Q\right)  }B\right)  +\underbrace{\sum_{\substack{R\subseteq
\left\{  1,2,\ldots,n\right\}  ;\\\left\vert R\right\vert =k;\\R\neq P}}0\det
A\cdot\det\left(  \operatorname*{sub}\nolimits_{w\left(  R\right)  }^{w\left(
Q\right)  }B\right)  }_{=0}\\
&  =\det A\cdot\det\left(  \operatorname*{sub}\nolimits_{w\left(  P\right)
}^{w\left(  Q\right)  }B\right)  .
\end{align*}
Hence, (\ref{pf.prop.sol.addexe.jacobi-complement.genform.11}) becomes%
\begin{align*}
&  \sum_{\substack{K\subseteq\left\{  1,2,\ldots,n\right\}  ;\\\left\vert
K\right\vert =\left\vert P\right\vert }}\left(  -1\right)  ^{\sum P+\sum
K}\det\left(  \operatorname*{sub}\nolimits_{w\left(  \widetilde{K}\right)
}^{w\left(  \widetilde{P}\right)  }A\right)  \det\left(  \operatorname*{sub}%
\nolimits_{w\left(  K\right)  }^{w\left(  Q\right)  }\left(  AB\right)
\right) \\
&  =\sum_{\substack{R\subseteq\left\{  1,2,\ldots,n\right\}  ;\\\left\vert
R\right\vert =k}}\delta_{R,P}\det A\cdot\det\left(  \operatorname*{sub}%
\nolimits_{w\left(  R\right)  }^{w\left(  Q\right)  }B\right)  =\det
A\cdot\det\left(  \operatorname*{sub}\nolimits_{w\left(  P\right)  }^{w\left(
Q\right)  }B\right)  .
\end{align*}
This proves Proposition \ref{prop.sol.addexe.jacobi-complement.genform}.
\end{proof}
\end{verlong}

Proposition \ref{prop.sol.addexe.jacobi-complement.genform} becomes
particularly simple when the matrix $AB$ is diagonal (i.e., has all entries
outside of its diagonal equal to $0$):

\begin{corollary}
\label{cor.sol.addexe.jacobi-complement.diagform}Let $n\in\mathbb{N}$. For
every two objects $i$ and $j$, define $\delta_{i,j}\in\mathbb{K}$ by
$\delta_{i,j}=%
\begin{cases}
1, & \text{if }i=j;\\
0, & \text{if }i\neq j
\end{cases}
$.

For any subset $I$ of $\left\{  1,2,\ldots,n\right\}  $, we let $\widetilde{I}%
$ denote the complement $\left\{  1,2,\ldots,n\right\}  \setminus I$ of $I$.

Let $d_{1},d_{2},\ldots,d_{n}$ be $n$ elements of $\mathbb{K}$. Let $A$ and
$B$ be two $n\times n$-matrices such that $AB=\left(  d_{i}\delta
_{i,j}\right)  _{1\leq i\leq n,\ 1\leq j\leq n}$.

Let $P$ be a subset of $\left\{  1,2,\ldots,n\right\}  $. Let $Q$ be a subset
of $\left\{  1,2,\ldots,n\right\}  $ such that $\left\vert P\right\vert
=\left\vert Q\right\vert $. Then,%
\[
\det A\cdot\det\left(  \operatorname*{sub}\nolimits_{w\left(  P\right)
}^{w\left(  Q\right)  }B\right)  =\left(  -1\right)  ^{\sum P+\sum Q}%
\det\left(  \operatorname*{sub}\nolimits_{w\left(  \widetilde{Q}\right)
}^{w\left(  \widetilde{P}\right)  }A\right)  \prod_{i\in Q}d_{i}.
\]

\end{corollary}

\begin{proof}
[Proof of Corollary \ref{cor.sol.addexe.jacobi-complement.diagform}.]If $K$ is
a subset of $\left\{  1,2,\ldots,n\right\}  $ satisfying $\left\vert
K\right\vert =\left\vert P\right\vert $, then%
\begin{equation}
\det\left(  \operatorname*{sub}\nolimits_{w\left(  K\right)  }^{w\left(
Q\right)  }\left(  AB\right)  \right)  =\delta_{K,Q}\prod_{i\in K}d_{i}
\label{pf.cor.sol.addexe.jacobi-complement.diagform.2}%
\end{equation}
\footnote{\textit{Proof of
(\ref{pf.cor.sol.addexe.jacobi-complement.diagform.2}):} Let $K$ be a subset
of $\left\{  1,2,\ldots,n\right\}  $ satisfying $\left\vert K\right\vert
=\left\vert P\right\vert $. Then, $\left\vert K\right\vert =\left\vert
P\right\vert =\left\vert Q\right\vert $. Also, $AB$ is the $n\times n$-matrix
$\left(  d_{i}\delta_{i,j}\right)  _{1\leq i\leq n,\ 1\leq j\leq n}$ (since
$AB=\left(  d_{i}\delta_{i,j}\right)  _{1\leq i\leq n,\ 1\leq j\leq n}$).
Hence, Lemma \ref{lem.diag.minors} (applied to $AB$ and $K$ instead of $D$ and
$P$) yields%
\[
\det\left(  \operatorname*{sub}\nolimits_{w\left(  K\right)  }^{w\left(
Q\right)  }\left(  AB\right)  \right)  =\delta_{K,Q}\prod_{i\in K}d_{i}.
\]
This proves (\ref{pf.cor.sol.addexe.jacobi-complement.diagform.2}).}.

Proposition \ref{prop.sol.addexe.jacobi-complement.genform} (applied to $m=n$)
yields%
\begin{align}
&  \det A\cdot\det\left(  \operatorname*{sub}\nolimits_{w\left(  P\right)
}^{w\left(  Q\right)  }B\right) \nonumber\\
&  =\sum_{\substack{K\subseteq\left\{  1,2,\ldots,n\right\}  ;\\\left\vert
K\right\vert =\left\vert P\right\vert }}\left(  -1\right)  ^{\sum P+\sum
K}\det\left(  \operatorname*{sub}\nolimits_{w\left(  \widetilde{K}\right)
}^{w\left(  \widetilde{P}\right)  }A\right)  \underbrace{\det\left(
\operatorname*{sub}\nolimits_{w\left(  K\right)  }^{w\left(  Q\right)
}\left(  AB\right)  \right)  }_{\substack{=\delta_{K,Q}\prod_{i\in K}%
d_{i}\\\text{(by (\ref{pf.cor.sol.addexe.jacobi-complement.diagform.2}))}%
}}\nonumber\\
&  =\sum_{\substack{K\subseteq\left\{  1,2,\ldots,n\right\}  ;\\\left\vert
K\right\vert =\left\vert P\right\vert }}\left(  -1\right)  ^{\sum P+\sum
K}\det\left(  \operatorname*{sub}\nolimits_{w\left(  \widetilde{K}\right)
}^{w\left(  \widetilde{P}\right)  }A\right)  \delta_{K,Q}\prod_{i\in K}d_{i}.
\label{pf.cor.sol.addexe.jacobi-complement.diagform.1}%
\end{align}

But $Q$ is a subset of $\left\{  1,2,\ldots,n\right\}  $ and satisfies
$\left\vert Q\right\vert =\left\vert P\right\vert $ (since $\left\vert
P\right\vert =\left\vert Q\right\vert $). In other words, $Q$ is a subset $K$
of $\left\{  1,2,\ldots,n\right\}  $ satisfying $\left\vert K\right\vert
=\left\vert P\right\vert $. Thus, the sum
\[
\sum_{\substack{K\subseteq\left\{  1,2,\ldots,n\right\}  ;\\\left\vert
K\right\vert =\left\vert P\right\vert }}\left(  -1\right)  ^{\sum P+\sum
K}\det\left(  \operatorname*{sub}\nolimits_{w\left(  \widetilde{K}\right)
}^{w\left(  \widetilde{P}\right)  }A\right)  \delta_{K,Q}\prod_{i\in K}d_{i}%
\]
has an addend for $K=Q$. If we split off this addend from this sum, then we
obtain%
\begin{align*}
&  \sum_{\substack{K\subseteq\left\{  1,2,\ldots,n\right\}  ;\\\left\vert
K\right\vert =\left\vert P\right\vert }}\left(  -1\right)  ^{\sum P+\sum
K}\det\left(  \operatorname*{sub}\nolimits_{w\left(  \widetilde{K}\right)
}^{w\left(  \widetilde{P}\right)  }A\right)  \delta_{K,Q}\prod_{i\in K}d_{i}\\
&  =\left(  -1\right)  ^{\sum P+\sum Q}\det\left(  \operatorname*{sub}%
\nolimits_{w\left(  \widetilde{Q}\right)  }^{w\left(  \widetilde{P}\right)
}A\right)  \underbrace{\delta_{Q,Q}}_{\substack{=1\\\text{(since }Q=Q\text{)}%
}}\prod_{i\in Q}d_{i}\\
&  \ \ \ \ \ \ \ \ \ \ +\sum_{\substack{K\subseteq\left\{  1,2,\ldots
,n\right\}  ;\\\left\vert K\right\vert =\left\vert P\right\vert ;\\K\neq
Q}}\left(  -1\right)  ^{\sum P+\sum K}\det\left(  \operatorname*{sub}%
\nolimits_{w\left(  \widetilde{K}\right)  }^{w\left(  \widetilde{P}\right)
}A\right)  \underbrace{\delta_{K,Q}}_{\substack{=0\\\text{(since }K\neq
Q\text{)}}}\prod_{i\in K}d_{i}\\
&  =\left(  -1\right)  ^{\sum P+\sum Q}\det\left(  \operatorname*{sub}%
\nolimits_{w\left(  \widetilde{Q}\right)  }^{w\left(  \widetilde{P}\right)
}A\right)  \prod_{i\in Q}d_{i}\\
&  \ \ \ \ \ \ \ \ \ \ +\underbrace{\sum_{\substack{K\subseteq\left\{
1,2,\ldots,n\right\}  ;\\\left\vert K\right\vert =\left\vert P\right\vert
;\\K\neq Q}}\left(  -1\right)  ^{\sum P+\sum K}\det\left(  \operatorname*{sub}%
\nolimits_{w\left(  \widetilde{K}\right)  }^{w\left(  \widetilde{P}\right)
}A\right)  0\prod_{i\in K}d_{i}}_{=0}\\
&  =\left(  -1\right)  ^{\sum P+\sum Q}\det\left(  \operatorname*{sub}%
\nolimits_{w\left(  \widetilde{Q}\right)  }^{w\left(  \widetilde{P}\right)
}A\right)  \prod_{i\in Q}d_{i}.
\end{align*}
Hence, (\ref{pf.cor.sol.addexe.jacobi-complement.diagform.1}) becomes%
\begin{align*}
&  \det A\cdot\det\left(  \operatorname*{sub}\nolimits_{w\left(  P\right)
}^{w\left(  Q\right)  }B\right) \\
&  =\sum_{\substack{K\subseteq\left\{  1,2,\ldots,n\right\}  ;\\\left\vert
K\right\vert =\left\vert P\right\vert }}\left(  -1\right)  ^{\sum P+\sum
K}\det\left(  \operatorname*{sub}\nolimits_{w\left(  \widetilde{K}\right)
}^{w\left(  \widetilde{P}\right)  }A\right)  \delta_{K,Q}\prod_{i\in K}d_{i}\\
&  =\left(  -1\right)  ^{\sum P+\sum Q}\det\left(  \operatorname*{sub}%
\nolimits_{w\left(  \widetilde{Q}\right)  }^{w\left(  \widetilde{P}\right)
}A\right)  \prod_{i\in Q}d_{i}.
\end{align*}
This proves Corollary \ref{cor.sol.addexe.jacobi-complement.diagform}.
\end{proof}

Specializing Corollary \ref{cor.sol.addexe.jacobi-complement.diagform} a bit
further, we obtain the following:

\begin{corollary}
\label{cor.sol.addexe.jacobi-complement.lidform}Let $n\in\mathbb{N}$. Let
$\lambda\in\mathbb{K}$.

For any subset $I$ of $\left\{  1,2,\ldots,n\right\}  $, we let $\widetilde{I}%
$ denote the complement $\left\{  1,2,\ldots,n\right\}  \setminus I$ of $I$.

Let $A$ and $B$ be two $n\times n$-matrices such that $AB=\lambda I_{n}$.

Let $P$ be a subset of $\left\{  1,2,\ldots,n\right\}  $. Let $Q$ be a subset
of $\left\{  1,2,\ldots,n\right\}  $ such that $\left\vert P\right\vert
=\left\vert Q\right\vert $. Then,%
\[
\det A\cdot\det\left(  \operatorname*{sub}\nolimits_{w\left(  P\right)
}^{w\left(  Q\right)  }B\right)  =\left(  -1\right)  ^{\sum P+\sum Q}%
\lambda^{\left\vert Q\right\vert }\det\left(  \operatorname*{sub}%
\nolimits_{w\left(  \widetilde{Q}\right)  }^{w\left(  \widetilde{P}\right)
}A\right)  .
\]

\end{corollary}

\begin{vershort}
\begin{proof}
[Proof of Corollary \ref{cor.sol.addexe.jacobi-complement.lidform}.]For every
two objects $i$ and $j$, define $\delta_{i,j}\in\mathbb{K}$ as in Corollary
\ref{cor.sol.addexe.jacobi-complement.diagform}. We have $I_{n}=\left(
\delta_{i,j}\right)  _{1\leq i\leq n,\ 1\leq j\leq n}$ (by the definition of
$I_{n}$). Now,%
\[
AB=\lambda\underbrace{I_{n}}_{\substack{=\left(  \delta_{i,j}\right)  _{1\leq
i\leq n,\ 1\leq j\leq n}}}=\lambda\left(  \delta_{i,j}\right)  _{1\leq i\leq
n,\ 1\leq j\leq n}=\left(  \lambda\delta_{i,j}\right)  _{1\leq i\leq n,\ 1\leq
j\leq n}.
\]
Hence, Corollary \ref{cor.sol.addexe.jacobi-complement.diagform} (applied to
$d_{k}=\lambda$) yields%
\begin{align*}
\det A\cdot\det\left(  \operatorname*{sub}\nolimits_{w\left(  P\right)
}^{w\left(  Q\right)  }B\right)   &  =\left(  -1\right)  ^{\sum P+\sum Q}%
\det\left(  \operatorname*{sub}\nolimits_{w\left(  \widetilde{Q}\right)
}^{w\left(  \widetilde{P}\right)  }A\right)  \underbrace{\prod_{i\in Q}%
\lambda}_{=\lambda^{\left\vert Q\right\vert }}\\
&  =\left(  -1\right)  ^{\sum P+\sum Q}\lambda^{\left\vert Q\right\vert }%
\det\left(  \operatorname*{sub}\nolimits_{w\left(  \widetilde{Q}\right)
}^{w\left(  \widetilde{P}\right)  }A\right)  .
\end{align*}
This proves Corollary \ref{cor.sol.addexe.jacobi-complement.lidform}.
\end{proof}
\end{vershort}

\begin{verlong}
\begin{proof}
[Proof of Corollary \ref{cor.sol.addexe.jacobi-complement.lidform}.]For every
two objects $i$ and $j$, define $\delta_{i,j}\in\mathbb{K}$ by $\delta_{i,j}=%
\begin{cases}
1, & \text{if }i=j;\\
0, & \text{if }i\neq j
\end{cases}
$. We have%
\[
AB=\lambda\underbrace{I_{n}}_{\substack{=\left(  \delta_{i,j}\right)  _{1\leq
i\leq n,\ 1\leq j\leq n}\\\text{(by the definition}\\\text{of }I_{n}\text{)}%
}}=\lambda\left(  \delta_{i,j}\right)  _{1\leq i\leq n,\ 1\leq j\leq
n}=\left(  \lambda\delta_{i,j}\right)  _{1\leq i\leq n,\ 1\leq j\leq n}.
\]
Hence, Corollary \ref{cor.sol.addexe.jacobi-complement.diagform} (applied to
$d_{k}=\lambda$) yields%
\begin{align*}
\det A\cdot\det\left(  \operatorname*{sub}\nolimits_{w\left(  P\right)
}^{w\left(  Q\right)  }B\right)   &  =\left(  -1\right)  ^{\sum P+\sum Q}%
\det\left(  \operatorname*{sub}\nolimits_{w\left(  \widetilde{Q}\right)
}^{w\left(  \widetilde{P}\right)  }A\right)  \underbrace{\prod_{i\in Q}%
\lambda}_{=\lambda^{\left\vert Q\right\vert }}\\
&  =\left(  -1\right)  ^{\sum P+\sum Q}\det\left(  \operatorname*{sub}%
\nolimits_{w\left(  \widetilde{Q}\right)  }^{w\left(  \widetilde{P}\right)
}A\right)  \lambda^{\left\vert Q\right\vert }\\
&  =\left(  -1\right)  ^{\sum P+\sum Q}\lambda^{\left\vert Q\right\vert }%
\det\left(  \operatorname*{sub}\nolimits_{w\left(  \widetilde{Q}\right)
}^{w\left(  \widetilde{P}\right)  }A\right)  .
\end{align*}
This proves Corollary \ref{cor.sol.addexe.jacobi-complement.lidform}.
\end{proof}
\end{verlong}

Finally, we obtain Exercise \ref{addexe.jacobi-complement} easily by setting
$\lambda=1$ in Corollary \ref{cor.sol.addexe.jacobi-complement.lidform}:

\begin{proof}
[Second solution to Exercise \ref{addexe.jacobi-complement}.]The matrix
$A\in\mathbb{K}^{n\times n}$ is invertible. Hence, the matrix $A^{-1}%
\in\mathbb{K}^{n\times n}$ is well-defined. Theorem \ref{thm.det(AB)} (applied
to $B=A^{-1}$) yields $\det\left(  AA^{-1}\right)  =\det A\cdot\det\left(
A^{-1}\right)  $. Thus,
\[
\det A\cdot\det\left(  A^{-1}\right)  =\det\left(  \underbrace{AA^{-1}%
}_{=I_{n}}\right)  =\det\left(  I_{n}\right)  =1.
\]

We have $A^{-1}A=I_{n}=1\cdot I_{n}$. Hence, Corollary
\ref{cor.sol.addexe.jacobi-complement.lidform} (applied to $A^{-1}$, $A$ and
$1$ instead of $A$, $B$ and $\lambda$) yields%
\begin{align*}
\det\left(  A^{-1}\right)  \cdot\det\left(  \operatorname*{sub}%
\nolimits_{w\left(  P\right)  }^{w\left(  Q\right)  }A\right)   &  =\left(
-1\right)  ^{\sum P+\sum Q}\underbrace{1^{\left\vert Q\right\vert }}_{=1}%
\det\left(  \operatorname*{sub}\nolimits_{w\left(  \widetilde{Q}\right)
}^{w\left(  \widetilde{P}\right)  }\left(  A^{-1}\right)  \right) \\
&  =\left(  -1\right)  ^{\sum P+\sum Q}\det\left(  \operatorname*{sub}%
\nolimits_{w\left(  \widetilde{Q}\right)  }^{w\left(  \widetilde{P}\right)
}\left(  A^{-1}\right)  \right)  .
\end{align*}
Multiplying both sides of this equality by $\det A$, we obtain%
\[
\det A\cdot\det\left(  A^{-1}\right)  \cdot\det\left(  \operatorname*{sub}%
\nolimits_{w\left(  P\right)  }^{w\left(  Q\right)  }A\right)  =\det
A\cdot\left(  -1\right)  ^{\sum P+\sum Q}\det\left(  \operatorname*{sub}%
\nolimits_{w\left(  \widetilde{Q}\right)  }^{w\left(  \widetilde{P}\right)
}\left(  A^{-1}\right)  \right)  .
\]
Comparing this with%
\[
\underbrace{\det A\cdot\det\left(  A^{-1}\right)  }_{=1}\cdot\det\left(
\operatorname*{sub}\nolimits_{w\left(  P\right)  }^{w\left(  Q\right)
}A\right)  =\det\left(  \operatorname*{sub}\nolimits_{w\left(  P\right)
}^{w\left(  Q\right)  }A\right)  ,
\]
we obtain%
\begin{align*}
\det\left(  \operatorname*{sub}\nolimits_{w\left(  P\right)  }^{w\left(
Q\right)  }A\right)   &  =\det A\cdot\left(  -1\right)  ^{\sum P+\sum Q}%
\det\left(  \operatorname*{sub}\nolimits_{w\left(  \widetilde{Q}\right)
}^{w\left(  \widetilde{P}\right)  }\left(  A^{-1}\right)  \right) \\
&  =\left(  -1\right)  ^{\sum P+\sum Q}\det A\cdot\det\left(
\operatorname*{sub}\nolimits_{w\left(  \widetilde{Q}\right)  }^{w\left(
\widetilde{P}\right)  }\left(  A^{-1}\right)  \right)  .
\end{align*}
Thus, Exercise \ref{addexe.jacobi-complement} is solved again.
\end{proof}

\subsubsection{Addendum}

As an easy consequence of our Second solution to Exercise
\ref{addexe.jacobi-complement}, we can obtain the following fact:

\begin{corollary}
\label{cor.sol.addexe.jacobi-complement.adjform}Let $n\in\mathbb{N}$. For any
subset $I$ of $\left\{  1,2,\ldots,n\right\}  $, we let $\widetilde{I}$ denote
the complement $\left\{  1,2,\ldots,n\right\}  \setminus I$ of $I$.

Let $A$ be an $n\times n$-matrix.

Let $P$ and $Q$ be two subsets of $\left\{  1,2,\ldots,n\right\}  $ such that
$\left\vert P\right\vert =\left\vert Q\right\vert $. Then,%
\[
\det A\cdot\det\left(  \operatorname*{sub}\nolimits_{w\left(  P\right)
}^{w\left(  Q\right)  }\left(  \operatorname*{adj}A\right)  \right)  =\left(
-1\right)  ^{\sum P+\sum Q}\left(  \det A\right)  ^{\left\vert Q\right\vert
}\det\left(  \operatorname*{sub}\nolimits_{w\left(  \widetilde{Q}\right)
}^{w\left(  \widetilde{P}\right)  }A\right)  .
\]

\end{corollary}

\begin{proof}
[Proof of Corollary \ref{cor.sol.addexe.jacobi-complement.adjform}.]Theorem
\ref{thm.adj.inverse} yields $A\cdot\operatorname*{adj}A=\operatorname*{adj}%
A\cdot A=\det A\cdot I_{n}$. Hence, Corollary
\ref{cor.sol.addexe.jacobi-complement.lidform} (applied to
$B=\operatorname*{adj}A$ and $\lambda=\det A$) yields%
\[
\det A\cdot\det\left(  \operatorname*{sub}\nolimits_{w\left(  P\right)
}^{w\left(  Q\right)  }\left(  \operatorname*{adj}A\right)  \right)  =\left(
-1\right)  ^{\sum P+\sum Q}\left(  \det A\right)  ^{\left\vert Q\right\vert
}\det\left(  \operatorname*{sub}\nolimits_{w\left(  \widetilde{Q}\right)
}^{w\left(  \widetilde{P}\right)  }A\right)  .
\]
This proves Corollary \ref{cor.sol.addexe.jacobi-complement.adjform}.
\end{proof}

Note that we could have also deduced Corollary
\ref{cor.sol.addexe.jacobi-complement.adjform} from the First solution to
Exercise \ref{addexe.jacobi-complement} (but this would have been more
difficult, since we would have to slightly generalize our argument).

It is possible to strengthen Corollary
\ref{cor.sol.addexe.jacobi-complement.adjform} in the case when $\left\vert
P\right\vert =\left\vert Q\right\vert \geq1$ as follows:

\begin{corollary}
\label{cor.sol.addexe.jacobi-complement.adjform.stronger}Let $n\in\mathbb{N}$.
For any subset $I$ of $\left\{  1,2,\ldots,n\right\}  $, we let $\widetilde{I}%
$ denote the complement $\left\{  1,2,\ldots,n\right\}  \setminus I$ of $I$.

Let $A$ be an $n\times n$-matrix.

Let $P$ and $Q$ be two subsets of $\left\{  1,2,\ldots,n\right\}  $ such that
$\left\vert P\right\vert =\left\vert Q\right\vert \geq1$. Then,%
\[
\det\left(  \operatorname*{sub}\nolimits_{w\left(  P\right)  }^{w\left(
Q\right)  }\left(  \operatorname*{adj}A\right)  \right)  =\left(  -1\right)
^{\sum P+\sum Q}\left(  \det A\right)  ^{\left\vert Q\right\vert -1}%
\det\left(  \operatorname*{sub}\nolimits_{w\left(  \widetilde{Q}\right)
}^{w\left(  \widetilde{P}\right)  }A\right)  .
\]

\end{corollary}

Loosely speaking, the claim of Corollary
\ref{cor.sol.addexe.jacobi-complement.adjform.stronger} is obtained from that
of Corollary \ref{cor.sol.addexe.jacobi-complement.adjform} by cancelling
$\det A$. However, it is not immediately clear that this cancellation is
allowed (for instance, $\det A$ could be $0$, or could be a nonzero
non-cancellable element). There are ways to justify this cancellation in full
generality; however, these are not in the scope of these notes.

\subsection{Solution to Exercise \ref{exe.det.pluecker.multi}}

Throughout this section, we shall use the notations introduced in Definition
\ref{def.submatrix}, in Definition \ref{def.unrows} and in Definition
\ref{def.sect.laplace.notations}. Also, whenever $m$ is an integer, we shall
use the notation $\left[  m\right]  $ for the set $\left\{  1,2,\ldots
,m\right\}  $. Furthermore, we shall use the following notations:

\begin{definition}
\label{def.sol.det.pluecker.multi.1}Let $n\in\mathbb{N}$ and $m\in\mathbb{N}$.
Let $B\in\mathbb{K}^{n\times m}$. Let $I$ be a subset of $\left\{
1,2,\ldots,m\right\}  $.

\textbf{(a)} Let $B_{\bullet,\sim I}$ denote the $n\times\left(  m-\left\vert
I\right\vert \right)  $-matrix whose columns are $B_{\bullet,j_{1}}%
,B_{\bullet,j_{2}},\ldots,B_{\bullet,j_{h}}$ (from left to right), where
$\left(  j_{1},j_{2},\ldots,j_{h}\right)  =w\left(  \left\{  1,2,\ldots
,m\right\}  \setminus I\right)  $. (We will see that this is well-defined in
Lemma \ref{lem.sol.det.pluecker.multi.1} \textbf{(a)} below.)

\textbf{(b)} Let $p\in\mathbb{N}$. Let $A\in\mathbb{K}^{n\times p}$. Let
$\left(  A\mid B_{\bullet,I}\right)  $ denote the $n\times\left(  p+\left\vert
I\right\vert \right)  $-matrix whose columns are $\underbrace{A_{\bullet
,1},A_{\bullet,2},\ldots,A_{\bullet,p}}_{\text{the columns of }A}%
,B_{\bullet,i_{1}},B_{\bullet,i_{2}},\ldots,B_{\bullet,i_{\ell}}$ (from left
to right), where $\left(  i_{1},i_{2},\ldots,i_{\ell}\right)  =w\left(
I\right)  $. (We will see that this is well-defined in Lemma
\ref{lem.sol.det.pluecker.multi.2} \textbf{(a)} below.)
\end{definition}

This definition agrees with the notations defined in Exercise
\ref{exe.det.pluecker.multi}, for the following reasons:

\begin{itemize}
\item The notation $B_{\bullet,\sim I}$ introduced in Definition
\ref{def.sol.det.pluecker.multi.1} \textbf{(a)} generalizes the notation
$B_{\bullet,\sim I}$ defined in Exercise \ref{exe.det.pluecker.multi}.
(Indeed, the former becomes the latter when we set $m=n+k$.)

\item The notation $\left(  A\mid B_{\bullet,I}\right)  $ introduced in
Definition \ref{def.sol.det.pluecker.multi.1} \textbf{(b)} generalizes the
notation $\left(  A\mid B_{\bullet,I}\right)  $ defined in Exercise
\ref{exe.det.pluecker.multi}. (Indeed, the former becomes the latter when we
set $m=n+k$ and $p=n-k$.)
\end{itemize}

We shall now state and prove two lemmas which show that the matrices
$B_{\bullet,\sim I}$ and $\left(  A\mid B_{\bullet,I}\right)  $ from
Definition \ref{def.sol.det.pluecker.multi.1} are well-defined, and express
these matrices in more familiar terms (which will be helpful when we come to
the solution of Exercise \ref{exe.det.pluecker.multi}). Both lemmas are fairly
obvious when you draw the matrices; but the rigorous proofs are rather tedious.

\begin{lemma}
\label{lem.sol.det.pluecker.multi.1}Let $n\in\mathbb{N}$ and $m\in\mathbb{N}$.
Let $B\in\mathbb{K}^{n\times m}$. Let $I$ be a subset of $\left\{
1,2,\ldots,m\right\}  $. Then:

\textbf{(a)} The matrix $B_{\bullet,\sim I}$ in Definition
\ref{def.sol.det.pluecker.multi.1} \textbf{(a)} is a well-defined
$n\times\left(  m-\left\vert I\right\vert \right)  $-matrix.

\textbf{(b)} We have
\[
B_{\bullet,\sim I}=\operatorname*{sub}\nolimits_{w\left(  \left[  n\right]
\right)  }^{w\left(  \left[  m\right]  \setminus I\right)  }B.
\]

\end{lemma}

(To make sense of the right-hand side of this equality, we recall that
$\left[  n\right]  =\left\{  1,2,\ldots,n\right\}  $ and $\left[  m\right]
=\left\{  1,2,\ldots,m\right\}  $, and that $\operatorname*{sub}%
\nolimits_{w\left(  \left[  n\right]  \right)  }^{w\left(  \left[  m\right]
\setminus I\right)  }B$ is defined as in Definition
\ref{def.sect.laplace.notations}.)

\begin{vershort}
\begin{proof}
[Proof of Lemma \ref{lem.sol.det.pluecker.multi.1}.]From $I\subseteq\left\{
1,2,\ldots,m\right\}  $, we obtain
\[
\left\vert \left\{  1,2,\ldots,m\right\}  \setminus I\right\vert
=\underbrace{\left\vert \left\{  1,2,\ldots,m\right\}  \right\vert }%
_{=m}-\left\vert I\right\vert =m-\left\vert I\right\vert .
\]

Define a list $\left(  j_{1},j_{2},\ldots,j_{h}\right)  $ by $\left(
j_{1},j_{2},\ldots,j_{h}\right)  =w\left(  \left\{  1,2,\ldots,m\right\}
\setminus I\right)  $. Then, the definition of $B_{\bullet,\sim I}$ says that
$B_{\bullet,\sim I}$ is the $n\times\left(  m-\left\vert I\right\vert \right)
$-matrix whose columns are $B_{\bullet,j_{1}},B_{\bullet,j_{2}},\ldots
,B_{\bullet,j_{h}}$ (from left to right). In order to check that this is
well-defined, we must verify that there really exists such an $n\times\left(
m-\left\vert I\right\vert \right)  $-matrix.

We have $\left(  j_{1},j_{2},\ldots,j_{h}\right)  =w\left(  \left\{
1,2,\ldots,m\right\}  \setminus I\right)  $. In other words, $\left(
j_{1},j_{2},\ldots,j_{h}\right)  $ is the list of all elements of $\left\{
1,2,\ldots,m\right\}  \setminus I$ in increasing order (with no repetitions)
(because this is how $w\left(  \left\{  1,2,\ldots,m\right\}  \setminus
I\right)  $ is defined). Hence, the length $h$ of this list equals $\left\vert
\left\{  1,2,\ldots,m\right\}  \setminus I\right\vert $. Thus, $h=\left\vert
\left\{  1,2,\ldots,m\right\}  \setminus I\right\vert =m-\left\vert
I\right\vert $.

Recall that $\left(  j_{1},j_{2},\ldots,j_{h}\right)  =w\left(  \left\{
1,2,\ldots,m\right\}  \setminus I\right)  $. Hence, $j_{1},j_{2},\ldots,j_{h}$
are $h$ elements of $\left\{  1,2,\ldots,m\right\}  \setminus I$, and thus are
elements of $\left\{  1,2,\ldots,m\right\}  $. Hence, $B_{\bullet,j_{1}%
},B_{\bullet,j_{2}},\ldots,B_{\bullet,j_{h}}$ are $h$ column vectors with $n$
entries each. Hence, there exists an $n\times h$-matrix whose columns are
$B_{\bullet,j_{1}},B_{\bullet,j_{2}},\ldots,B_{\bullet,j_{h}}$ (from left to
right). In view of $h=m-\left\vert I\right\vert $, this rewrites as follows:
There exists an $n\times\left(  m-\left\vert I\right\vert \right)  $-matrix
whose columns are $B_{\bullet,j_{1}},B_{\bullet,j_{2}},\ldots,B_{\bullet
,j_{h}}$ (from left to right). In other words, the matrix $B_{\bullet,\sim I}$
in Definition \ref{def.sol.det.pluecker.multi.1} \textbf{(a)} is a
well-defined $n\times\left(  m-\left\vert I\right\vert \right)  $-matrix. This
proves Lemma \ref{lem.sol.det.pluecker.multi.1} \textbf{(a)}.

\textbf{(b)} Write the $n\times m$-matrix $B$ in the form $B=\left(
b_{i,j}\right)  _{1\leq i\leq n,\ 1\leq j\leq m}$. The list $w\left(  \left[
n\right]  \right)  $ is defined as the list of all elements of $\left[
n\right]  $ in increasing order (with no repetitions), and thus is simply the
list $\left(  1,2,\ldots,n\right)  $; in other words, $w\left(  \left[
n\right]  \right)  =\left(  1,2,\ldots,n\right)  $. Moreover, $w\left(
\underbrace{\left[  m\right]  }_{=\left\{  1,2,\ldots,m\right\}  }\setminus
I\right)  =w\left(  \left\{  1,2,\ldots,m\right\}  \setminus I\right)
=\left(  j_{1},j_{2},\ldots,j_{h}\right)  $. Now,%
\begin{align}
\operatorname*{sub}\nolimits_{w\left(  \left[  n\right]  \right)  }^{w\left(
\left[  m\right]  \setminus I\right)  }B  &  =\operatorname*{sub}%
\nolimits_{\left(  1,2,\ldots,n\right)  }^{\left(  j_{1},j_{2},\ldots
,j_{h}\right)  }B\ \ \ \ \ \ \ \ \ \ \left(
\begin{array}
[c]{c}%
\text{since }w\left(  \left[  n\right]  \right)  =\left(  1,2,\ldots,n\right)
\\
\text{and }w\left(  \left[  m\right]  \setminus I\right)  =\left(  j_{1}%
,j_{2},\ldots,j_{h}\right)
\end{array}
\right) \nonumber\\
&  =\operatorname*{sub}\nolimits_{1,2,\ldots,n}^{j_{1},j_{2},\ldots,j_{h}%
}B=\left(  b_{x,j_{y}}\right)  _{1\leq x\leq n,\ 1\leq y\leq h}
\label{pf.lem.sol.det.pluecker.multi.1.b.short.1}%
\end{align}
(by the definition of $\operatorname*{sub}\nolimits_{1,2,\ldots,n}%
^{j_{1},j_{2},\ldots,j_{h}}B$, since $B=\left(  b_{i,j}\right)  _{1\leq i\leq
n,\ 1\leq j\leq m}$).

On the other hand, the definition of $B_{\bullet,\sim I}$ shows that the
columns of the matrix $B_{\bullet,\sim I}$ are $B_{\bullet,j_{1}}%
,B_{\bullet,j_{2}},\ldots,B_{\bullet,j_{h}}$ (from left to right). Thus, for
each $y\in\left\{  1,2,\ldots,h\right\}  $, we have%
\begin{equation}
\left(  \text{the }y\text{-th column of the matrix }B_{\bullet,\sim I}\right)
=B_{\bullet,j_{y}}. \label{pf.lem.sol.det.pluecker.multi.1.b.short.q-col}%
\end{equation}
In particular, the matrix $B_{\bullet,\sim I}$ has $h$ columns, each of which
is a column vector with $n$ entries; thus, $B_{\bullet,\sim I}$ is an $n\times
h$-matrix. For each $x\in\left\{  1,2,\ldots,n\right\}  $ and $y\in\left\{
1,2,\ldots,h\right\}  $, we have
\begin{align*}
&  \left(  \text{the }\left(  x,y\right)  \text{-th entry of the matrix
}B_{\bullet,\sim I}\right) \\
&  =\left(  \text{the }x\text{-th entry of the }y\text{-th column of the
matrix }B_{\bullet,\sim I}\right) \\
&  =\left(  \text{the }x\text{-th entry of }B_{\bullet,j_{y}}\right)
\ \ \ \ \ \ \ \ \ \ \left(  \text{by
(\ref{pf.lem.sol.det.pluecker.multi.1.b.short.q-col})}\right) \\
&  =\left(  \text{the }\left(  x,j_{y}\right)  \text{-th entry of }B\right)
\ \ \ \ \ \ \ \ \ \ \left(  \text{since }B_{\bullet,j_{y}}\text{ is the }%
j_{y}\text{-th column of }B\right) \\
&  =b_{x,j_{y}}\ \ \ \ \ \ \ \ \ \ \left(  \text{since }B=\left(
b_{i,j}\right)  _{1\leq i\leq n,\ 1\leq j\leq m}\right)  .
\end{align*}
Hence, $B_{\bullet,\sim I}=\left(  b_{x,j_{y}}\right)  _{1\leq x\leq n,\ 1\leq
y\leq h}$. Comparing this with
(\ref{pf.lem.sol.det.pluecker.multi.1.b.short.1}), we obtain $B_{\bullet,\sim
I}=\operatorname*{sub}\nolimits_{w\left(  \left[  n\right]  \right)
}^{w\left(  \left[  m\right]  \setminus I\right)  }B$. This proves Lemma
\ref{lem.sol.det.pluecker.multi.1} \textbf{(b)}.
\end{proof}
\end{vershort}

\begin{verlong}
\begin{proof}
[Proof of Lemma \ref{lem.sol.det.pluecker.multi.1}.]We have $\left[  m\right]
=\left\{  1,2,\ldots,m\right\}  $ (by the definition of $\left[  m\right]  $).
From $I\subseteq\left\{  1,2,\ldots,m\right\}  =\left[  m\right]  $, we obtain
$\left\vert \left[  m\right]  \setminus I\right\vert =\underbrace{\left\vert
\left[  m\right]  \right\vert }_{=m}-\left\vert I\right\vert =m-\left\vert
I\right\vert $.

Also, $\left[  n\right]  =\left\{  1,2,\ldots,n\right\}  $ (by the definition
of $\left[  n\right]  $).

Define a list $\left(  j_{1},j_{2},\ldots,j_{h}\right)  $ by $\left(
j_{1},j_{2},\ldots,j_{h}\right)  =w\left(  \left\{  1,2,\ldots,m\right\}
\setminus I\right)  $. Then, the definition of $B_{\bullet,\sim I}$ says that
$B_{\bullet,\sim I}$ is the $n\times\left(  m-\left\vert I\right\vert \right)
$-matrix whose columns are $B_{\bullet,j_{1}},B_{\bullet,j_{2}},\ldots
,B_{\bullet,j_{h}}$ (from left to right). In order to check that this is
well-defined, we must verify that there really exists an $n\times\left(
m-\left\vert I\right\vert \right)  $-matrix whose columns are $B_{\bullet
,j_{1}},B_{\bullet,j_{2}},\ldots,B_{\bullet,j_{h}}$ (from left to right).

First of all,
\[
\left(  j_{1},j_{2},\ldots,j_{h}\right)  =w\left(  \underbrace{\left\{
1,2,\ldots,m\right\}  }_{=\left[  m\right]  }\setminus I\right)  =w\left(
\left[  m\right]  \setminus I\right)  .
\]

Proposition \ref{prop.sect.laplace.notations.w(I)} (applied to $\left[
m\right]  \setminus I$ instead of $I$) shows that $w\left(  \left[  m\right]
\setminus I\right)  $ is an $\left\vert \left[  m\right]  \setminus
I\right\vert $-tuple of elements of $\left[  m\right]  \setminus I$. In other
words, $\left(  j_{1},j_{2},\ldots,j_{h}\right)  $ is an $\left\vert \left[
m\right]  \setminus I\right\vert $-tuple of elements of $\left[  m\right]
\setminus I$ (since $\left(  j_{1},j_{2},\ldots,j_{h}\right)  =w\left(
\left[  m\right]  \setminus I\right)  $). Thus, $\left(  j_{1},j_{2}%
,\ldots,j_{h}\right)  $ is an $\left\vert \left[  m\right]  \setminus
I\right\vert $-tuple, i.e., a list of length $\left\vert \left[  m\right]
\setminus I\right\vert $. Hence, the length of the list $\left(  j_{1}%
,j_{2},\ldots,j_{h}\right)  $ is $\left\vert \left[  m\right]  \setminus
I\right\vert $. Thus,
\[
\left(  \text{the length of the list }\left(  j_{1},j_{2},\ldots,j_{h}\right)
\right)  =\left\vert \left[  m\right]  \setminus I\right\vert =m-\left\vert
I\right\vert .
\]
Comparing this with%
\[
\left(  \text{the length of the list }\left(  j_{1},j_{2},\ldots,j_{h}\right)
\right)  =h,
\]
we obtain $h=m-\left\vert I\right\vert $.

Now, $B_{\bullet,j_{p}}$ is a well-defined column vector with $n$ entries for
each $p\in\left\{  1,2,\ldots,h\right\}  $\ \ \ \ \footnote{\textit{Proof.}
Let $p\in\left\{  1,2,\ldots,h\right\}  $. Then, $j_{p}$ is an element of
$\left[  m\right]  \setminus I$ (since $\left(  j_{1},j_{2},\ldots
,j_{h}\right)  $ is an $\left\vert \left[  m\right]  \setminus I\right\vert
$-tuple of elements of $\left[  m\right]  \setminus I$). Hence, $j_{p}%
\in\left[  m\right]  \setminus I\subseteq\left[  m\right]  =\left\{
1,2,\ldots,m\right\}  $. Thus, $B_{\bullet,j_{p}}$ is a well-defined column
vector with $n$ entries (since $B$ is an $n\times m$-matrix). Qed.}. In other
words, $B_{\bullet,j_{1}},B_{\bullet,j_{2}},\ldots,B_{\bullet,j_{h}}$ are $h$
column vectors with $n$ entries each. Hence, there exists an $n\times
h$-matrix whose columns are $B_{\bullet,j_{1}},B_{\bullet,j_{2}}%
,\ldots,B_{\bullet,j_{h}}$ (from left to right). In view of $h=m-\left\vert
I\right\vert $, this rewrites as follows: There exists an $n\times\left(
m-\left\vert I\right\vert \right)  $-matrix whose columns are $B_{\bullet
,j_{1}},B_{\bullet,j_{2}},\ldots,B_{\bullet,j_{h}}$ (from left to right). In
other words, the matrix $B_{\bullet,\sim I}$ in Definition
\ref{def.sol.det.pluecker.multi.1} \textbf{(a)} is a well-defined
$n\times\left(  m-\left\vert I\right\vert \right)  $-matrix. This proves Lemma
\ref{lem.sol.det.pluecker.multi.1} \textbf{(a)}.

\textbf{(b)} The set $\left[  m\right]  \setminus I$ is a subset of $\left\{
1,2,\ldots,m\right\}  $ (since $\left[  m\right]  \setminus I\subseteq\left[
m\right]  =\left\{  1,2,\ldots,m\right\}  $), whereas the set $\left[
n\right]  $ is a subset of $\left\{  1,2,\ldots,n\right\}  $ (since $\left[
n\right]  =\left\{  1,2,\ldots,n\right\}  $). Hence, the matrix
$\operatorname*{sub}\nolimits_{w\left(  \left[  n\right]  \right)  }^{w\left(
\left[  m\right]  \setminus I\right)  }B$ is well-defined.

The definition of $B_{\bullet,\sim I}$ says that $B_{\bullet,\sim I}$ is the
$n\times\left(  m-\left\vert I\right\vert \right)  $-matrix whose columns are
$B_{\bullet,j_{1}},B_{\bullet,j_{2}},\ldots,B_{\bullet,j_{h}}$ (from left to
right). Thus, the columns of the matrix $B_{\bullet,\sim I}$ are
$B_{\bullet,j_{1}},B_{\bullet,j_{2}},\ldots,B_{\bullet,j_{h}}$ (from left to
right). Hence, for each $q\in\left\{  1,2,\ldots,h\right\}  $, the $q$-th
column of the matrix $B_{\bullet,\sim I}$ is $B_{\bullet,j_{q}}$. In other
words, for each $q\in\left\{  1,2,\ldots,h\right\}  $, we have%
\begin{equation}
\left(  \text{the }q\text{-th column of the matrix }B_{\bullet,\sim I}\right)
=B_{\bullet,j_{q}}. \label{pf.lem.sol.det.pluecker.multi.1.col1}%
\end{equation}

We have $w\left(  \left[  n\right]  \right)  =\left(  1,2,\ldots,n\right)
$\ \ \ \ \footnote{\textit{Proof.} Recall that $w\left(  \left[  n\right]
\right)  $ is the list of all elements of $\left[  n\right]  $ in increasing
order (with no repetitions) (by the definition of $w\left(  \left[  n\right]
\right)  $). Thus,%
\begin{align*}
w\left(  \left[  n\right]  \right)   &  =\left(  \text{the list of all
elements of }\underbrace{\left[  n\right]  }_{=\left\{  1,2,\ldots,n\right\}
}\text{ in increasing order (with no repetitions)}\right) \\
&  =\left(  \text{the list of all elements of }\left\{  1,2,\ldots,n\right\}
\text{ in increasing order (with no repetitions)}\right) \\
&  =\left(  1,2,\ldots,n\right)  .
\end{align*}
}. Now, write the $n\times m$-matrix $B$ in the form $B=\left(  b_{i,j}%
\right)  _{1\leq i\leq n,\ 1\leq j\leq m}$. Then,%
\begin{align*}
\operatorname*{sub}\nolimits_{w\left(  \left[  n\right]  \right)  }^{w\left(
\left[  m\right]  \setminus I\right)  }B  &  =\operatorname*{sub}%
\nolimits_{\left(  1,2,\ldots,n\right)  }^{\left(  j_{1},j_{2},\ldots
,j_{h}\right)  }B\ \ \ \ \ \ \ \ \ \ \left(
\begin{array}
[c]{c}%
\text{since }w\left(  \left[  n\right]  \right)  =\left(  1,2,\ldots,n\right)
\\
\text{and }w\left(  \left[  m\right]  \setminus I\right)  =\left(  j_{1}%
,j_{2},\ldots,j_{h}\right)
\end{array}
\right) \\
&  =\operatorname*{sub}\nolimits_{1,2,\ldots,n}^{j_{1},j_{2},\ldots,j_{h}%
}B=\left(  b_{x,j_{y}}\right)  _{1\leq x\leq n,\ 1\leq y\leq h}\\
&  \ \ \ \ \ \ \ \ \ \ \left(  \text{by the definition of }\operatorname*{sub}%
\nolimits_{1,2,\ldots,n}^{j_{1},j_{2},\ldots,j_{h}}B\text{, since }B=\left(
b_{i,j}\right)  _{1\leq i\leq n,\ 1\leq j\leq m}\right)  .
\end{align*}
Hence, $\operatorname*{sub}\nolimits_{w\left(  \left[  n\right]  \right)
}^{w\left(  \left[  m\right]  \setminus I\right)  }B$ is an $n\times
h$-matrix. Also, $B_{\bullet,\sim I}$ is an $n\times\left(  m-\left\vert
I\right\vert \right)  $-matrix; in other words, $B_{\bullet,\sim I}$ is an
$n\times h$-matrix (since $h=m-\left\vert I\right\vert $).

Let $q\in\left\{  1,2,\ldots,h\right\}  $. Then, $\left(  B_{\bullet,\sim
I}\right)  _{\bullet,q}$ and $\left(  \operatorname*{sub}\nolimits_{w\left(
\left[  n\right]  \right)  }^{w\left(  \left[  m\right]  \setminus I\right)
}B\right)  _{\bullet,q}$ are well-defined (since both $B_{\bullet,\sim I}$ and
$\operatorname*{sub}\nolimits_{w\left(  \left[  n\right]  \right)  }^{w\left(
\left[  m\right]  \setminus I\right)  }B$ are $n\times h$-matrices). Moreover,
$\left(  B_{\bullet,\sim I}\right)  _{\bullet,q}$ is the $q$-th column of the
matrix $B_{\bullet,\sim I}$ (by the definition of $\left(  B_{\bullet,\sim
I}\right)  _{\bullet,q}$). Hence,%
\begin{align}
\left(  B_{\bullet,\sim I}\right)  _{\bullet,q}  &  =\left(  \text{the
}q\text{-th column of the matrix }B_{\bullet,\sim I}\right)  =B_{\bullet
,j_{q}}\ \ \ \ \ \ \ \ \ \ \left(  \text{by
(\ref{pf.lem.sol.det.pluecker.multi.1.col1})}\right) \nonumber\\
&  =\left(  \text{the }j_{q}\text{-th column of the matrix }\underbrace{B}%
_{=\left(  b_{i,j}\right)  _{1\leq i\leq n,\ 1\leq j\leq m}}\right)
\nonumber\\
&  \ \ \ \ \ \ \ \ \ \ \left(  \text{since }B_{\bullet,j_{q}}\text{ is defined
as the }j_{q}\text{-th column of the matrix }B\right) \nonumber\\
&  =\left(  \text{the }j_{q}\text{-th column of the matrix }\left(
b_{i,j}\right)  _{1\leq i\leq n,\ 1\leq j\leq m}\right) \nonumber\\
&  =\left(
\begin{array}
[c]{c}%
b_{1,j_{q}}\\
b_{2,j_{q}}\\
\vdots\\
b_{n,j_{q}}%
\end{array}
\right)  . \label{pf.lem.sol.det.pluecker.multi.1.col2}%
\end{align}

On the other hand, $\operatorname*{sub}\nolimits_{w\left(  \left[  n\right]
\right)  }^{w\left(  \left[  m\right]  \setminus I\right)  }B$ is an $n\times
h$-matrix. Hence, $\left(  \operatorname*{sub}\nolimits_{w\left(  \left[
n\right]  \right)  }^{w\left(  \left[  m\right]  \setminus I\right)
}B\right)  _{\bullet,q}$ is the $q$-th column of the matrix
$\operatorname*{sub}\nolimits_{w\left(  \left[  n\right]  \right)  }^{w\left(
\left[  m\right]  \setminus I\right)  }B$ (by the definition of $\left(
\operatorname*{sub}\nolimits_{w\left(  \left[  n\right]  \right)  }^{w\left(
\left[  m\right]  \setminus I\right)  }B\right)  _{\bullet,q}$). Thus,%
\begin{align*}
\left(  \operatorname*{sub}\nolimits_{w\left(  \left[  n\right]  \right)
}^{w\left(  \left[  m\right]  \setminus I\right)  }B\right)  _{\bullet,q}  &
=\left(  \text{the }q\text{-th column of the matrix }%
\underbrace{\operatorname*{sub}\nolimits_{w\left(  \left[  n\right]  \right)
}^{w\left(  \left[  m\right]  \setminus I\right)  }B}_{=\left(  b_{x,j_{y}%
}\right)  _{1\leq x\leq n,\ 1\leq y\leq h}}\right) \\
&  =\left(  \text{the }q\text{-th column of the matrix }\left(  b_{x,j_{y}%
}\right)  _{1\leq x\leq n,\ 1\leq y\leq h}\right) \\
&  =\left(
\begin{array}
[c]{c}%
b_{1,j_{q}}\\
b_{2,j_{q}}\\
\vdots\\
b_{n,j_{q}}%
\end{array}
\right)  .
\end{align*}
Comparing this with (\ref{pf.lem.sol.det.pluecker.multi.1.col2}), we obtain
$\left(  B_{\bullet,\sim I}\right)  _{\bullet,q}=\left(  \operatorname*{sub}%
\nolimits_{w\left(  \left[  n\right]  \right)  }^{w\left(  \left[  m\right]
\setminus I\right)  }B\right)  _{\bullet,q}$.

Now, forget that we fixed $q$. We thus have shown that $\left(  B_{\bullet
,\sim I}\right)  _{\bullet,q}=\left(  \operatorname*{sub}\nolimits_{w\left(
\left[  n\right]  \right)  }^{w\left(  \left[  m\right]  \setminus I\right)
}B\right)  _{\bullet,q}$ for each $q\in\left\{  1,2,\ldots,h\right\}  $. Thus,
Lemma \ref{lem.sol.prop.addcol.props.cols} (applied to $h$, $B_{\bullet,\sim
I}$ and $\operatorname*{sub}\nolimits_{w\left(  \left[  n\right]  \right)
}^{w\left(  \left[  m\right]  \setminus I\right)  }B$ instead of $m$, $A$ and
$B$) shows that $B_{\bullet,\sim I}=\operatorname*{sub}\nolimits_{w\left(
\left[  n\right]  \right)  }^{w\left(  \left[  m\right]  \setminus I\right)
}B$ (since both $B_{\bullet,\sim I}$ and $\operatorname*{sub}%
\nolimits_{w\left(  \left[  n\right]  \right)  }^{w\left(  \left[  m\right]
\setminus I\right)  }B$ are $n\times h$-matrices). This proves Lemma
\ref{lem.sol.det.pluecker.multi.1} \textbf{(b)}.
\end{proof}
\end{verlong}

\begin{lemma}
\label{lem.sol.det.pluecker.multi.2}Let $n\in\mathbb{N}$ and $m\in\mathbb{N}$
and $p\in\mathbb{N}$. Let $B\in\mathbb{K}^{n\times m}$. Let $I$ be a subset of
$\left\{  1,2,\ldots,m\right\}  $. Let $A\in\mathbb{K}^{n\times p}$.

\textbf{(a)} The matrix $\left(  A\mid B_{\bullet,I}\right)  $ in Definition
\ref{def.sol.det.pluecker.multi.1} \textbf{(b)} is a well-defined
$n\times\left(  p+\left\vert I\right\vert \right)  $-matrix.

\textbf{(b)} Let $P$ be any subset of $\left\{  1,2,\ldots,n\right\}  $. Let
$Q$ be the set $\left\{  1,2,\ldots,p\right\}  $. Let $R$ be the set $\left\{
p+1,p+2,\ldots,p+\left\vert I\right\vert \right\}  $. Then,%
\[
\operatorname*{sub}\nolimits_{w\left(  P\right)  }^{w\left(  Q\right)
}\left(  A\mid B_{\bullet,I}\right)  =\operatorname*{sub}\nolimits_{w\left(
P\right)  }^{w\left(  Q\right)  }A
\]
and%
\[
\operatorname*{sub}\nolimits_{w\left(  P\right)  }^{w\left(  R\right)
}\left(  A\mid B_{\bullet,I}\right)  =\operatorname*{sub}\nolimits_{w\left(
P\right)  }^{w\left(  I\right)  }B.
\]

\end{lemma}

\begin{vershort}
\begin{proof}
[Proof of Lemma \ref{lem.sol.det.pluecker.multi.2}.]Define a list $\left(
i_{1},i_{2},\ldots,i_{\ell}\right)  $ by $\left(  i_{1},i_{2},\ldots,i_{\ell
}\right)  =w\left(  I\right)  $. Then, the definition of $\left(  A\mid
B_{\bullet,I}\right)  $ says that $\left(  A\mid B_{\bullet,I}\right)  $ is
the $n\times\left(  p+\left\vert I\right\vert \right)  $-matrix whose columns
are $\underbrace{A_{\bullet,1},A_{\bullet,2},\ldots,A_{\bullet,p}}_{\text{the
columns of }A},B_{\bullet,i_{1}},B_{\bullet,i_{2}},\ldots,B_{\bullet,i_{\ell}%
}$ (from left to right). In order to check that this is well-defined, we must
verify that there really exists such an $n\times\left(  p+\left\vert
I\right\vert \right)  $-matrix.

We know that $A\in\mathbb{K}^{n\times p}$. Hence, $A_{\bullet,1},A_{\bullet
,2},\ldots,A_{\bullet,p}$ are $p$ column vectors in $\mathbb{K}^{n\times1}$:
namely, the columns of $A$.

We have $\left(  i_{1},i_{2},\ldots,i_{\ell}\right)  =w\left(  I\right)  $. In
other words, $\left(  i_{1},i_{2},\ldots,i_{\ell}\right)  $ is the list of all
elements of $I$ in increasing order (with no repetition) (because this is how
$w\left(  I\right)  $ is defined). Thus, the length $\ell$ of this list equals
$\left\vert I\right\vert $. In other words, $\ell=\left\vert I\right\vert $.

Also, $i_{1},i_{2},\ldots,i_{\ell}$ are elements of $I$ (since $\left(
i_{1},i_{2},\ldots,i_{\ell}\right)  $ is the list of all elements of $I$ in
increasing order), and thus are elements of $\left\{  1,2,\ldots,m\right\}  $
(since $I$ is a subset of $\left\{  1,2,\ldots,m\right\}  $). Hence,
$B_{\bullet,i_{1}},B_{\bullet,i_{2}},\ldots,B_{\bullet,i_{\ell}}$ are $\ell$
column vectors in $\mathbb{K}^{n\times1}$ (since $B\in\mathbb{K}^{n\times m}$).

We now conclude that $\underbrace{A_{\bullet,1},A_{\bullet,2},\ldots
,A_{\bullet,p}}_{\text{the columns of }A},B_{\bullet,i_{1}},B_{\bullet,i_{2}%
},\ldots,B_{\bullet,i_{\ell}}$ are $p+\ell$ column vectors in $\mathbb{K}%
^{n\times1}$ (because $A_{\bullet,1},A_{\bullet,2},\ldots,A_{\bullet,p}$ are
$p$ column vectors in $\mathbb{K}^{n\times1}$, and because $B_{\bullet,i_{1}%
},B_{\bullet,i_{2}},\ldots,B_{\bullet,i_{\ell}}$ are $\ell$ column vectors in
$\mathbb{K}^{n\times1}$). Thus, there exists an $n\times\left(  p+\ell\right)
$-matrix whose columns are $\underbrace{A_{\bullet,1},A_{\bullet,2}%
,\ldots,A_{\bullet,p}}_{\text{the columns of }A},B_{\bullet,i_{1}}%
,B_{\bullet,i_{2}},\ldots,B_{\bullet,i_{\ell}}$ (from left to right). In view
of $\ell=\left\vert I\right\vert $, this rewrites as follows: There exists an
$n\times\left(  p+\left\vert I\right\vert \right)  $-matrix whose columns are
$\underbrace{A_{\bullet,1},A_{\bullet,2},\ldots,A_{\bullet,p}}_{\text{the
columns of }A},B_{\bullet,i_{1}},B_{\bullet,i_{2}},\ldots,B_{\bullet,i_{\ell}%
}$ (from left to right). In other words, the matrix $\left(  A\mid
B_{\bullet,I}\right)  $ in Definition \ref{def.sol.det.pluecker.multi.1}
\textbf{(b)} is a well-defined $n\times\left(  p+\left\vert I\right\vert
\right)  $-matrix. This proves Lemma \ref{lem.sol.det.pluecker.multi.2}
\textbf{(a)}.

\textbf{(b)} The definition of $R$ yields $R=\left\{  p+1,p+2,\ldots
,p+\left\vert I\right\vert \right\}  \subseteq\left\{  1,2,\ldots,p+\left\vert
I\right\vert \right\}  $ (since $p\geq0$).

The definition of $Q$ yields $Q=\left\{  1,2,\ldots,p\right\}  \subseteq
\left\{  1,2,\ldots,p+\left\vert I\right\vert \right\}  $ (since $\left\vert
I\right\vert \geq0$).

Recall that $\left(  A\mid B_{\bullet,I}\right)  $ is an $n\times\left(
p+\left\vert I\right\vert \right)  $-matrix. In other words, $\left(  A\mid
B_{\bullet,I}\right)  $ is an $n\times\left(  p+\ell\right)  $-matrix (since
$\ell=\left\vert I\right\vert $). Write this $n\times\left(  p+\ell\right)
$-matrix $\left(  A\mid B_{\bullet,I}\right)  $ in the form $\left(  A\mid
B_{\bullet,I}\right)  =\left(  c_{i,j}\right)  _{1\leq i\leq n,\ 1\leq j\leq
p+\ell}$.

Write the $n\times p$-matrix $A$ in the form $A=\left(  a_{i,j}\right)
_{1\leq i\leq n,\ 1\leq j\leq p}$.

Write the $n\times m$-matrix $B$ in the form $B=\left(  b_{i,j}\right)
_{1\leq i\leq n,\ 1\leq j\leq m}$.

The columns of the matrix $\left(  A\mid B_{\bullet,I}\right)  $ (from left to
right) are \newline$A_{\bullet,1},A_{\bullet,2},\ldots,A_{\bullet
,p},B_{\bullet,i_{1}},B_{\bullet,i_{2}},\ldots,B_{\bullet,i_{\ell}}$ (by the
definition of $\left(  A\mid B_{\bullet,I}\right)  $). Thus, for each
$g\in\left\{  1,2,\ldots,p+\ell\right\}  $, we have%
\begin{equation}
\left(  \text{the }g\text{-th column of the matrix }\left(  A\mid
B_{\bullet,I}\right)  \right)  =%
\begin{cases}
A_{\bullet,g}, & \text{if }g\leq p;\\
B_{\bullet,i_{g-p}}, & \text{if }g>p
\end{cases}
. \label{pf.lem.sol.det.pluecker.multi.2.short.b.gcol}%
\end{equation}

For each $x\in\left\{  1,2,\ldots,n\right\}  $ and $g\in\left\{
1,2,\ldots,p\right\}  $, we have%
\begin{equation}
c_{x,g}=a_{x,g} \label{pf.lem.sol.det.pluecker.multi.2.short.b.cxg}%
\end{equation}
\footnote{\textit{Proof of (\ref{pf.lem.sol.det.pluecker.multi.2.short.b.cxg}%
):} Let $x\in\left\{  1,2,\ldots,n\right\}  $ and $g\in\left\{  1,2,\ldots
,p\right\}  $. Thus, $g\leq p$. Also, $g\in\left\{  1,2,\ldots,p\right\}
\subseteq\left\{  1,2,\ldots,p+\ell\right\}  $ (since $\ell\geq0$). Hence,
(\ref{pf.lem.sol.det.pluecker.multi.2.short.b.gcol}) yields%
\begin{align}
&  \left(  \text{the }g\text{-th column of the matrix }\left(  A\mid
B_{\bullet,I}\right)  \right) \nonumber\\
&  =%
\begin{cases}
A_{\bullet,g}, & \text{if }g\leq p;\\
B_{\bullet,i_{g-p}}, & \text{if }g>p
\end{cases}
=A_{\bullet,g}\ \ \ \ \ \ \ \ \ \ \left(  \text{since }g\leq p\right)
\nonumber\\
&  =\left(  \text{the }g\text{-th column of the matrix }A\right)
\label{pf.lem.sol.det.pluecker.multi.2.short.b.cxg.pf.1}%
\end{align}
(since $A_{\bullet,g}$ is defined to be the $g$-th column of the matrix $A$).
But
\[
\left(  \text{the }\left(  x,g\right)  \text{-th entry of the matrix }\left(
A\mid B_{\bullet,I}\right)  \right)  =c_{x,g}%
\]
(since $\left(  A\mid B_{\bullet,I}\right)  =\left(  c_{i,j}\right)  _{1\leq
i\leq n,\ 1\leq j\leq p+\ell}$). Comparing this with%
\begin{align*}
&  \left(  \text{the }\left(  x,g\right)  \text{-th entry of the matrix
}\left(  A\mid B_{\bullet,I}\right)  \right) \\
&  =\left(  \text{the }x\text{-th entry of }\underbrace{\text{the }g\text{-th
column of the matrix }\left(  A\mid B_{\bullet,I}\right)  }%
_{\substack{=\left(  \text{the }g\text{-th column of the matrix }A\right)
\\\text{(by (\ref{pf.lem.sol.det.pluecker.multi.2.short.b.cxg.pf.1}))}%
}}\right) \\
&  =\left(  \text{the }x\text{-th entry of the }g\text{-th column of the
matrix }A\right) \\
&  =\left(  \text{the }\left(  x,g\right)  \text{-th entry of the matrix
}A\right)  =a_{x,g}\ \ \ \ \ \ \ \ \ \ \left(  \text{since }A=\left(
a_{i,j}\right)  _{1\leq i\leq n,\ 1\leq j\leq p}\right)  ,
\end{align*}
we obtain $c_{x,g}=a_{x,g}$. This proves
(\ref{pf.lem.sol.det.pluecker.multi.2.short.b.cxg}).}.

For each $x\in\left\{  1,2,\ldots,n\right\}  $ and $g\in\left\{
p+1,p+2,\ldots,p+\ell\right\}  $, we have%
\begin{equation}
c_{x,g}=b_{x,i_{g-p}} \label{pf.lem.sol.det.pluecker.multi.2.short.b.cxg2}%
\end{equation}
\footnote{\textit{Proof of (\ref{pf.lem.sol.det.pluecker.multi.2.short.b.cxg2}%
):} Let $x\in\left\{  1,2,\ldots,n\right\}  $ and $g\in\left\{  p+1,p+2,\ldots
,p+\ell\right\}  $. Thus, $g\geq p+1>p$. Also, $g\in\left\{  p+1,p+2,\ldots
,p+\ell\right\}  \subseteq\left\{  1,2,\ldots,p+\ell\right\}  $. Hence,
(\ref{pf.lem.sol.det.pluecker.multi.2.short.b.gcol}) yields%
\begin{align}
&  \left(  \text{the }g\text{-th column of the matrix }\left(  A\mid
B_{\bullet,I}\right)  \right) \nonumber\\
&  =%
\begin{cases}
A_{\bullet,g}, & \text{if }g\leq p;\\
B_{\bullet,i_{g-p}}, & \text{if }g>p
\end{cases}
=B_{\bullet,i_{g-p}}\ \ \ \ \ \ \ \ \ \ \left(  \text{since }g>p\right)
\nonumber\\
&  =\left(  \text{the }i_{g-p}\text{-th column of the matrix }B\right)
\label{pf.lem.sol.det.pluecker.multi.2.short.b.cxg2.pf.1}%
\end{align}
(since $B_{\bullet,i_{g-p}}$ is defined to be the $i_{g-p}$-th column of the
matrix $B$). But
\[
\left(  \text{the }\left(  x,g\right)  \text{-th entry of the matrix }\left(
A\mid B_{\bullet,I}\right)  \right)  =c_{x,g}%
\]
(since $\left(  A\mid B_{\bullet,I}\right)  =\left(  c_{i,j}\right)  _{1\leq
i\leq n,\ 1\leq j\leq p+\ell}$). Comparing this with%
\begin{align*}
&  \left(  \text{the }\left(  x,g\right)  \text{-th entry of the matrix
}\left(  A\mid B_{\bullet,I}\right)  \right) \\
&  =\left(  \text{the }x\text{-th entry of }\underbrace{\text{the }g\text{-th
column of the matrix }\left(  A\mid B_{\bullet,I}\right)  }%
_{\substack{=\left(  \text{the }i_{g-p}\text{-th column of the matrix
}B\right)  \\\text{(by
(\ref{pf.lem.sol.det.pluecker.multi.2.short.b.cxg2.pf.1}))}}}\right) \\
&  =\left(  \text{the }x\text{-th entry of the }i_{g-p}\text{-th column of the
matrix }B\right) \\
&  =\left(  \text{the }\left(  x,i_{g-p}\right)  \text{-th entry of the matrix
}B\right)  =b_{x,i_{g-p}}\ \ \ \ \ \ \ \ \ \ \left(  \text{since }B=\left(
b_{i,j}\right)  _{1\leq i\leq n,\ 1\leq j\leq m}\right)  ,
\end{align*}
we obtain $c_{x,g}=b_{x,i_{g-p}}$. This proves
(\ref{pf.lem.sol.det.pluecker.multi.2.short.b.cxg2}).}. Hence, for each
$x\in\left\{  1,2,\ldots,n\right\}  $ and $y\in\left\{  1,2,\ldots
,\ell\right\}  $, we have%
\begin{equation}
c_{x,p+y}=b_{x,i_{y}} \label{pf.lem.sol.det.pluecker.multi.2.short.b.cxg2b}%
\end{equation}
\footnote{\textit{Proof of
(\ref{pf.lem.sol.det.pluecker.multi.2.short.b.cxg2b}):} Let $x\in\left\{
1,2,\ldots,n\right\}  $ and $y\in\left\{  1,2,\ldots,\ell\right\}  $. From
$y\in\left\{  1,2,\ldots,\ell\right\}  $, we obtain $p+y\in\left\{
p+1,p+2,\ldots,p+\ell\right\}  $. Hence,
(\ref{pf.lem.sol.det.pluecker.multi.2.short.b.cxg2}) (applied to $g=p+y$)
yields $c_{x,p+y}=b_{x,i_{p+y-p}}=b_{x,i_{y}}$ (since $p+y-p=y$). This proves
(\ref{pf.lem.sol.det.pluecker.multi.2.short.b.cxg2b}).}.

Now, $P\subseteq\left\{  1,2,\ldots,n\right\}  $, whereas $Q\subseteq\left\{
1,2,\ldots,p+\left\vert I\right\vert \right\}  $. Thus, the matrix
\newline$\operatorname*{sub}\nolimits_{w\left(  P\right)  }^{w\left(
Q\right)  }\left(  A\mid B_{\bullet,I}\right)  $ is well-defined (since
$\left(  A\mid B_{\bullet,I}\right)  $ is an $n\times\left(  p+\left\vert
I\right\vert \right)  $-matrix). Also, $P\subseteq$ $\left\{  1,2,\ldots
,n\right\}  $, whereas $Q\subseteq\left\{  1,2,\ldots,p\right\}  $ (since
$Q=\left\{  1,2,\ldots,p\right\}  $). Hence, the matrix $\operatorname*{sub}%
\nolimits_{w\left(  P\right)  }^{w\left(  Q\right)  }A$ is well-defined (since
$A$ is an $n\times p$-matrix).

Furthermore, $P\subseteq\left\{  1,2,\ldots,n\right\}  $, whereas
$R\subseteq\left\{  1,2,\ldots,p+\left\vert I\right\vert \right\}  $. Thus,
the matrix $\operatorname*{sub}\nolimits_{w\left(  P\right)  }^{w\left(
R\right)  }\left(  A\mid B_{\bullet,I}\right)  $ is well-defined (since
$\left(  A\mid B_{\bullet,I}\right)  $ is an $n\times\left(  p+\left\vert
I\right\vert \right)  $-matrix). Also, $P\subseteq\left\{  1,2,\ldots
,n\right\}  $, whereas $I\subseteq\left\{  1,2,\ldots,m\right\}  $. Hence, the
matrix $\operatorname*{sub}\nolimits_{w\left(  P\right)  }^{w\left(  I\right)
}B$ is well-defined (since $B$ is an $n\times m$-matrix).

We have $Q=\left\{  1,2,\ldots,p\right\}  $, and thus $w\left(  Q\right)
=\left(  1,2,\ldots,p\right)  $ (by the definition of $w\left(  Q\right)  $).

Also, $R=\left\{  p+1,p+2,\ldots,p+\left\vert I\right\vert \right\}  =\left\{
p+1,p+2,\ldots,p+\ell\right\}  $ (since $\left\vert I\right\vert =\ell$), and
thus $w\left(  R\right)  =\left(  p+1,p+2,\ldots,p+\ell\right)  $ (by the
definition of $w\left(  R\right)  $).

Write the list $w\left(  P\right)  $ in the form $w\left(  P\right)  =\left(
p_{1},p_{2},\ldots,p_{h}\right)  $. Thus, $p_{1},p_{2},\ldots,p_{h}$ are all
the elements of $P$ (since $w\left(  P\right)  $ is the list of all elements
of $P$ in increasing order). Hence, $p_{1},p_{2},\ldots,p_{h}$ are elements of
$P$, and therefore are elements of $\left\{  1,2,\ldots,n\right\}  $ (since
$P\subseteq\left\{  1,2,\ldots,n\right\}  $). In other words, $p_{x}%
\in\left\{  1,2,\ldots,n\right\}  $ for each $x\in\left\{  1,2,\ldots
,h\right\}  $. Therefore, for each $x\in\left\{  1,2,\ldots,h\right\}  $ and
$y\in\left\{  1,2,\ldots,p\right\}  $, we have $c_{p_{x},y}=a_{p_{x},y}$ (by
(\ref{pf.lem.sol.det.pluecker.multi.2.short.b.cxg}) (applied to $p_{x}$ and
$y$ instead of $x$ and $g$)). In other words, we have%
\begin{equation}
\left(  c_{p_{x},y}\right)  _{1\leq x\leq h,\ 1\leq y\leq p}=\left(
a_{p_{x},y}\right)  _{1\leq x\leq h,\ 1\leq y\leq p}.
\label{pf.lem.sol.det.pluecker.multi.2.short.b.cpxy2a}%
\end{equation}

Recall that $p_{x}\in\left\{  1,2,\ldots,n\right\}  $ for each $x\in\left\{
1,2,\ldots,h\right\}  $. Therefore, for each $x\in\left\{  1,2,\ldots
,h\right\}  $ and $y\in\left\{  1,2,\ldots,\ell\right\}  $, we have
$c_{p_{x},p+y}=b_{p_{x},i_{y}}$ (by
(\ref{pf.lem.sol.det.pluecker.multi.2.short.b.cxg2b}) (applied to $p_{x}$
instead of $x$)). In other words, we have%
\begin{equation}
\left(  c_{p_{x},p+y}\right)  _{1\leq x\leq h,\ 1\leq y\leq\ell}=\left(
b_{p_{x},i_{y}}\right)  _{1\leq x\leq h,\ 1\leq y\leq\ell}.
\label{pf.lem.sol.det.pluecker.multi.2.short.b.cpxy2b}%
\end{equation}

Using $w\left(  P\right)  =\left(  p_{1},p_{2},\ldots,p_{h}\right)  $ and
$w\left(  Q\right)  =\left(  1,2,\ldots,p\right)  $, we see that
\begin{align*}
\operatorname*{sub}\nolimits_{w\left(  P\right)  }^{w\left(  Q\right)
}\left(  A\mid B_{\bullet,I}\right)   &  =\operatorname*{sub}%
\nolimits_{\left(  p_{1},p_{2},\ldots,p_{h}\right)  }^{\left(  1,2,\ldots
,p\right)  }\left(  A\mid B_{\bullet,I}\right)  =\operatorname*{sub}%
\nolimits_{p_{1},p_{2},\ldots,p_{h}}^{1,2,\ldots,p}\left(  A\mid B_{\bullet
,I}\right) \\
&  =\left(  c_{p_{x},y}\right)  _{1\leq x\leq h,\ 1\leq y\leq p}\\
&  \ \ \ \ \ \ \ \ \ \ \left(
\begin{array}
[c]{c}%
\text{by the definition of }\operatorname*{sub}\nolimits_{p_{1},p_{2}%
,\ldots,p_{h}}^{1,2,\ldots,p}\left(  A\mid B_{\bullet,I}\right)  \text{,}\\
\text{since }\left(  A\mid B_{\bullet,I}\right)  =\left(  c_{i,j}\right)
_{1\leq i\leq n,\ 1\leq j\leq p+\ell}%
\end{array}
\right) \\
&  =\left(  a_{p_{x},y}\right)  _{1\leq x\leq h,\ 1\leq y\leq p}%
\ \ \ \ \ \ \ \ \ \ \left(  \text{by
(\ref{pf.lem.sol.det.pluecker.multi.2.short.b.cpxy2a})}\right)  .
\end{align*}
Comparing this with%
\begin{align*}
\operatorname*{sub}\nolimits_{w\left(  P\right)  }^{w\left(  Q\right)  }A  &
=\operatorname*{sub}\nolimits_{\left(  p_{1},p_{2},\ldots,p_{h}\right)
}^{\left(  1,2,\ldots,p\right)  }A\ \ \ \ \ \ \ \ \ \ \left(
\begin{array}
[c]{c}%
\text{since }w\left(  P\right)  =\left(  p_{1},p_{2},\ldots,p_{h}\right) \\
\text{and }w\left(  Q\right)  =\left(  1,2,\ldots,p\right)
\end{array}
\right) \\
&  =\operatorname*{sub}\nolimits_{p_{1},p_{2},\ldots,p_{h}}^{1,2,\ldots
,p}A=\left(  a_{p_{x},y}\right)  _{1\leq x\leq h,\ 1\leq y\leq p}\\
&  \ \ \ \ \ \ \ \ \ \ \left(
\begin{array}
[c]{c}%
\text{by the definition of }\operatorname*{sub}\nolimits_{p_{1},p_{2}%
,\ldots,p_{h}}^{1,2,\ldots,p}A\text{,}\\
\text{since }A=\left(  a_{i,j}\right)  _{1\leq i\leq n,\ 1\leq j\leq p}%
\end{array}
\right)  ,
\end{align*}
we obtain%
\[
\operatorname*{sub}\nolimits_{w\left(  P\right)  }^{w\left(  Q\right)
}\left(  A\mid B_{\bullet,I}\right)  =\operatorname*{sub}\nolimits_{w\left(
P\right)  }^{w\left(  Q\right)  }A.
\]

Using $w\left(  P\right)  =\left(  p_{1},p_{2},\ldots,p_{h}\right)  $ and
$w\left(  R\right)  =\left(  p+1,p+2,\ldots,p+\ell\right)  $, we find%
\begin{align*}
\operatorname*{sub}\nolimits_{w\left(  P\right)  }^{w\left(  R\right)
}\left(  A\mid B_{\bullet,I}\right)   &  =\operatorname*{sub}%
\nolimits_{\left(  p_{1},p_{2},\ldots,p_{h}\right)  }^{\left(  p+1,p+2,\ldots
,p+\ell\right)  }\left(  A\mid B_{\bullet,I}\right)  =\operatorname*{sub}%
\nolimits_{p_{1},p_{2},\ldots,p_{h}}^{p+1,p+2,\ldots,p+\ell}\left(  A\mid
B_{\bullet,I}\right) \\
&  =\left(  c_{p_{x},p+y}\right)  _{1\leq x\leq h,\ 1\leq y\leq\ell}\\
&  \ \ \ \ \ \ \ \ \ \ \left(
\begin{array}
[c]{c}%
\text{by the definition of }\operatorname*{sub}\nolimits_{p_{1},p_{2}%
,\ldots,p_{h}}^{p+1,p+2,\ldots,p+\ell}\left(  A\mid B_{\bullet,I}\right)
\text{,}\\
\text{since }\left(  A\mid B_{\bullet,I}\right)  =\left(  c_{i,j}\right)
_{1\leq i\leq n,\ 1\leq j\leq p+\ell}%
\end{array}
\right) \\
&  =\left(  b_{p_{x},i_{y}}\right)  _{1\leq x\leq h,\ 1\leq y\leq\ell
}\ \ \ \ \ \ \ \ \ \ \left(  \text{by
(\ref{pf.lem.sol.det.pluecker.multi.2.short.b.cpxy2b})}\right)  .
\end{align*}
Comparing this with%
\begin{align*}
\operatorname*{sub}\nolimits_{w\left(  P\right)  }^{w\left(  I\right)  }B  &
=\operatorname*{sub}\nolimits_{\left(  p_{1},p_{2},\ldots,p_{h}\right)
}^{\left(  i_{1},i_{2},\ldots,i_{\ell}\right)  }B\ \ \ \ \ \ \ \ \ \ \left(
\begin{array}
[c]{c}%
\text{since }w\left(  I\right)  =\left(  i_{1},i_{2},\ldots,i_{\ell}\right) \\
\text{and }w\left(  P\right)  =\left(  p_{1},p_{2},\ldots,p_{h}\right)
\end{array}
\right) \\
&  =\operatorname*{sub}\nolimits_{p_{1},p_{2},\ldots,p_{h}}^{i_{1}%
,i_{2},\ldots,i_{\ell}}B=\left(  b_{p_{x},i_{y}}\right)  _{1\leq x\leq
h,\ 1\leq y\leq\ell}\\
&  \ \ \ \ \ \ \ \ \ \ \left(
\begin{array}
[c]{c}%
\text{by the definition of }\operatorname*{sub}\nolimits_{p_{1},p_{2}%
,\ldots,p_{h}}^{i_{1},i_{2},\ldots,i_{\ell}}B\text{,}\\
\text{since }B=\left(  b_{i,j}\right)  _{1\leq i\leq n,\ 1\leq j\leq m}%
\end{array}
\right)  ,
\end{align*}
we obtain%
\[
\operatorname*{sub}\nolimits_{w\left(  P\right)  }^{w\left(  R\right)
}\left(  A\mid B_{\bullet,I}\right)  =\operatorname*{sub}\nolimits_{w\left(
P\right)  }^{w\left(  I\right)  }B.
\]
Thus, the proof of Lemma \ref{lem.sol.det.pluecker.multi.2} \textbf{(b)} is complete.
\end{proof}
\end{vershort}

\begin{verlong}
\begin{proof}
[Proof of Lemma \ref{lem.sol.det.pluecker.multi.2}.]Define a list $\left(
i_{1},i_{2},\ldots,i_{\ell}\right)  $ by $\left(  i_{1},i_{2},\ldots,i_{\ell
}\right)  =w\left(  I\right)  $. Then, the definition of $\left(  A\mid
B_{\bullet,I}\right)  $ says that $\left(  A\mid B_{\bullet,I}\right)  $ is
the $n\times\left(  p+\left\vert I\right\vert \right)  $-matrix whose columns
are $\underbrace{A_{\bullet,1},A_{\bullet,2},\ldots,A_{\bullet,p}}_{\text{the
columns of }A},B_{\bullet,i_{1}},B_{\bullet,i_{2}},\ldots,B_{\bullet,i_{\ell}%
}$ (from left to right). In order to check that this is well-defined, we must
verify that there really exists an $n\times\left(  p+\left\vert I\right\vert
\right)  $-matrix whose columns are $\underbrace{A_{\bullet,1},A_{\bullet
,2},\ldots,A_{\bullet,p}}_{\text{the columns of }A},B_{\bullet,i_{1}%
},B_{\bullet,i_{2}},\ldots,B_{\bullet,i_{\ell}}$ (from left to right).

We know that $A$ is an $n\times p$-matrix (since $A\in\mathbb{K}^{n\times p}%
$). Hence, $A_{\bullet,1},A_{\bullet,2},\ldots,A_{\bullet,p}$ are $p$
well-defined column vectors in $\mathbb{K}^{n\times1}$: namely, these $p$
column vectors are the columns of $A$.

Proposition \ref{prop.sect.laplace.notations.w(I)} shows that $w\left(
I\right)  $ is an $\left\vert I\right\vert $-tuple of elements of $I$. In
other words, $\left(  i_{1},i_{2},\ldots,i_{\ell}\right)  $ is an $\left\vert
I\right\vert $-tuple of elements of $I$ (since $\left(  i_{1},i_{2}%
,\ldots,i_{\ell}\right)  =w\left(  I\right)  $). Thus, $\left(  i_{1}%
,i_{2},\ldots,i_{\ell}\right)  $ is an $\left\vert I\right\vert $-tuple, i.e.,
a list of length $\left\vert I\right\vert $. Hence, the length of the list
$\left(  i_{1},i_{2},\ldots,i_{\ell}\right)  $ is $\left\vert I\right\vert $.
Thus,
\[
\left(  \text{the length of the list }\left(  i_{1},i_{2},\ldots,i_{\ell
}\right)  \right)  =\left\vert I\right\vert .
\]
Comparing this with%
\[
\left(  \text{the length of the list }\left(  i_{1},i_{2},\ldots,i_{\ell
}\right)  \right)  =\ell,
\]
we obtain $\ell=\left\vert I\right\vert $.

Now, $B_{\bullet,i_{q}}$ is a well-defined column vector in $\mathbb{K}%
^{n\times1}$ for each $q\in\left\{  1,2,\ldots,\ell\right\}  $%
\ \ \ \ \footnote{\textit{Proof.} Let $q\in\left\{  1,2,\ldots,\ell\right\}
$. Then, $i_{q}$ is an element of $I$ (since $\left(  i_{1},i_{2}%
,\ldots,i_{\ell}\right)  $ is an $\left\vert I\right\vert $-tuple of elements
of $I$). Hence, $i_{q}\in I\subseteq\left\{  1,2,\ldots,m\right\}  $. Thus,
$B_{\bullet,i_{q}}$ is a well-defined column vector with $n$ entries (since
$B$ is an $n\times m$-matrix). In other words, $B_{\bullet,i_{q}}$ is a
well-defined column vector in $\mathbb{K}^{n\times1}$. Qed.}. In other words,
$B_{\bullet,i_{1}},B_{\bullet,i_{2}},\ldots,B_{\bullet,i_{\ell}}$ are $\ell$
well-defined column vectors in $\mathbb{K}^{n\times1}$.

We now conclude that $\underbrace{A_{\bullet,1},A_{\bullet,2},\ldots
,A_{\bullet,p}}_{\text{the columns of }A},B_{\bullet,i_{1}},B_{\bullet,i_{2}%
},\ldots,B_{\bullet,i_{\ell}}$ are $p+\ell$ well-defined column vectors in
$\mathbb{K}^{n\times1}$ (because $A_{\bullet,1},A_{\bullet,2},\ldots
,A_{\bullet,p}$ are $p$ well-defined column vectors in $\mathbb{K}^{n\times1}%
$, and because $B_{\bullet,i_{1}},B_{\bullet,i_{2}},\ldots,B_{\bullet,i_{\ell
}}$ are $\ell$ well-defined column vectors in $\mathbb{K}^{n\times1}$). Thus,
there exists an $n\times\left(  p+\ell\right)  $-matrix whose columns
are\newline$\underbrace{A_{\bullet,1},A_{\bullet,2},\ldots,A_{\bullet,p}%
}_{\text{the columns of }A},B_{\bullet,i_{1}},B_{\bullet,i_{2}},\ldots
,B_{\bullet,i_{\ell}}$ (from left to right). In view of $\ell=\left\vert
I\right\vert $, this rewrites as follows: There exists an $n\times\left(
p+\left\vert I\right\vert \right)  $-matrix whose columns are\newline%
$\underbrace{A_{\bullet,1},A_{\bullet,2},\ldots,A_{\bullet,p}}_{\text{the
columns of }A},B_{\bullet,i_{1}},B_{\bullet,i_{2}},\ldots,B_{\bullet,i_{\ell}%
}$ (from left to right). In other words, the matrix $\left(  A\mid
B_{\bullet,I}\right)  $ in Definition \ref{def.sol.det.pluecker.multi.1}
\textbf{(b)} is a well-defined $n\times\left(  p+\left\vert I\right\vert
\right)  $-matrix. This proves Lemma \ref{lem.sol.det.pluecker.multi.2}
\textbf{(a)}.

\textbf{(b)} We have $p\in\mathbb{N}$, so that $p\geq0$ and thus
$p+1\geq0+1=1$. Now, the definition of $R$ yields $R=\left\{  p+1,p+2,\ldots
,p+\left\vert I\right\vert \right\}  \subseteq\left\{  1,2,\ldots,p+\left\vert
I\right\vert \right\}  $ (since $p+1\geq1$).

Also, $\left\vert I\right\vert \in\mathbb{N}$, so that $\left\vert
I\right\vert \geq0$ and thus $p+\left\vert I\right\vert \geq p+0=p$. Hence,
$p\leq p+\left\vert I\right\vert $. Now, the definition of $Q$ yields
$Q=\left\{  1,2,\ldots,p\right\}  \subseteq\left\{  1,2,\ldots,p+\left\vert
I\right\vert \right\}  $ (since $p\leq p+\left\vert I\right\vert $).

The columns of the matrix $\left(  A\mid B_{\bullet,I}\right)  $ (from left to
right) are \newline$\underbrace{A_{\bullet,1},A_{\bullet,2},\ldots
,A_{\bullet,p}}_{\text{the columns of }A},B_{\bullet,i_{1}},B_{\bullet,i_{2}%
},\ldots,B_{\bullet,i_{\ell}}$\ \ \ \ \footnote{because $\left(  A\mid
B_{\bullet,I}\right)  $ is the $n\times\left(  p+\left\vert I\right\vert
\right)  $-matrix whose columns are $\underbrace{A_{\bullet,1},A_{\bullet
,2},\ldots,A_{\bullet,p}}_{\text{the columns of }A},B_{\bullet,i_{1}%
},B_{\bullet,i_{2}},\ldots,B_{\bullet,i_{\ell}}$ (from left to right)}. In
other words, the columns of the matrix $\left(  A\mid B_{\bullet,I}\right)  $
(from left to right) are $A_{\bullet,1},A_{\bullet,2},\ldots,A_{\bullet
,p},B_{\bullet,i_{1}},B_{\bullet,i_{2}},\ldots,B_{\bullet,i_{\ell}}$. Also,
$\left(  A\mid B_{\bullet,I}\right)  $ is an $n\times\left(  p+\left\vert
I\right\vert \right)  $-matrix. In other words, $\left(  A\mid B_{\bullet
,I}\right)  $ is an $n\times\left(  p+\ell\right)  $-matrix (since
$\ell=\left\vert I\right\vert $). Write this $n\times\left(  p+\ell\right)
$-matrix $\left(  A\mid B_{\bullet,I}\right)  $ in the form $\left(  A\mid
B_{\bullet,I}\right)  =\left(  c_{i,j}\right)  _{1\leq i\leq n,\ 1\leq j\leq
p+\ell}$.

Now, for each $g\in\left\{  1,2,\ldots,p+\ell\right\}  $, we have%
\begin{align}
&  \left(  \text{the }g\text{-th column of the matrix }\left(  A\mid
B_{\bullet,I}\right)  \right) \nonumber\\
&  =\left(  \text{the }g\text{-th among the columns of the matrix }\left(
A\mid B_{\bullet,I}\right)  \text{ (from left to right)}\right) \nonumber\\
&  =\left(  \text{the }g\text{-th among the column vectors }A_{\bullet
,1},A_{\bullet,2},\ldots,A_{\bullet,p},B_{\bullet,i_{1}},B_{\bullet,i_{2}%
},\ldots,B_{\bullet,i_{\ell}}\right) \nonumber\\
&  \ \ \ \ \ \ \ \ \ \ \left(
\begin{array}
[c]{c}%
\text{since the columns of the matrix }\left(  A\mid B_{\bullet,I}\right)
\text{ (from left to right)}\\
\text{are the column vectors }A_{\bullet,1},A_{\bullet,2},\ldots,A_{\bullet
,p},B_{\bullet,i_{1}},B_{\bullet,i_{2}},\ldots,B_{\bullet,i_{\ell}}%
\end{array}
\right) \nonumber\\
&  =%
\begin{cases}
A_{\bullet,g}, & \text{if }g\leq p;\\
B_{\bullet,i_{g-p}}, & \text{if }g>p
\end{cases}
. \label{pf.lem.sol.det.pluecker.multi.2.b.gcol}%
\end{align}

Write the $n\times p$-matrix $A$ in the form $A=\left(  a_{i,j}\right)
_{1\leq i\leq n,\ 1\leq j\leq p}$.

Write the $n\times m$-matrix $B$ in the form $B=\left(  b_{i,j}\right)
_{1\leq i\leq n,\ 1\leq j\leq m}$.

For each $x\in\left\{  1,2,\ldots,n\right\}  $ and $g\in\left\{
1,2,\ldots,p\right\}  $, we have%
\begin{equation}
c_{x,g}=a_{x,g} \label{pf.lem.sol.det.pluecker.multi.2.b.cxg}%
\end{equation}
\footnote{\textit{Proof of (\ref{pf.lem.sol.det.pluecker.multi.2.b.cxg}):} Let
$x\in\left\{  1,2,\ldots,n\right\}  $ and $g\in\left\{  1,2,\ldots,p\right\}
$. Then, $g\in\left\{  1,2,\ldots,p\right\}  \subseteq\left\{  1,2,\ldots
,p+\ell\right\}  $ (since $p\leq p+\ell$ (because $\ell\geq0$)). Hence,
(\ref{pf.lem.sol.det.pluecker.multi.2.b.gcol}) yields%
\begin{align}
&  \left(  \text{the }g\text{-th column of the matrix }\left(  A\mid
B_{\bullet,I}\right)  \right) \nonumber\\
&  =%
\begin{cases}
A_{\bullet,g}, & \text{if }g\leq p;\\
B_{\bullet,i_{g-p}}, & \text{if }g>p
\end{cases}
=A_{\bullet,g}\ \ \ \ \ \ \ \ \ \ \left(  \text{since }g\leq p\text{ (since
}g\in\left\{  1,2,\ldots,p\right\}  \text{)}\right) \nonumber\\
&  =\left(  \text{the }g\text{-th column of the matrix }A\right)
\label{pf.lem.sol.det.pluecker.multi.2.b.cxg.pf.1}%
\end{align}
(since $A_{\bullet,g}$ is defined to be the $g$-th column of the matrix $A$).
But
\[
\left(  \text{the }\left(  x,g\right)  \text{-th entry of the matrix }\left(
A\mid B_{\bullet,I}\right)  \right)  =c_{x,g}%
\]
(since $\left(  A\mid B_{\bullet,I}\right)  =\left(  c_{i,j}\right)  _{1\leq
i\leq n,\ 1\leq j\leq p+\ell}$). Comparing this with%
\begin{align*}
&  \left(  \text{the }\left(  x,g\right)  \text{-th entry of the matrix
}\left(  A\mid B_{\bullet,I}\right)  \right) \\
&  =\left(  \text{the }x\text{-th entry of }\underbrace{\text{the }g\text{-th
column of the matrix }\left(  A\mid B_{\bullet,I}\right)  }%
_{\substack{=\left(  \text{the }g\text{-th column of the matrix }A\right)
\\\text{(by (\ref{pf.lem.sol.det.pluecker.multi.2.b.cxg.pf.1}))}}}\right) \\
&  =\left(  \text{the }x\text{-th entry of the }g\text{-th column of the
matrix }A\right) \\
&  =\left(  \text{the }\left(  x,g\right)  \text{-th entry of the matrix
}A\right)  =a_{x,g}\ \ \ \ \ \ \ \ \ \ \left(  \text{since }A=\left(
a_{i,j}\right)  _{1\leq i\leq n,\ 1\leq j\leq p}\right)  ,
\end{align*}
we obtain $c_{x,g}=a_{x,g}$. This proves
(\ref{pf.lem.sol.det.pluecker.multi.2.b.cxg}).}.

For each $x\in\left\{  1,2,\ldots,n\right\}  $ and $g\in\left\{
p+1,p+2,\ldots,p+\ell\right\}  $, we have%
\begin{equation}
c_{x,g}=b_{x,i_{g-p}} \label{pf.lem.sol.det.pluecker.multi.2.b.cxg2}%
\end{equation}
\footnote{\textit{Proof of (\ref{pf.lem.sol.det.pluecker.multi.2.b.cxg2}):}
Let $x\in\left\{  1,2,\ldots,n\right\}  $ and $g\in\left\{  p+1,p+2,\ldots
,p+\ell\right\}  $. Then, $g\in\left\{  p+1,p+2,\ldots,p+\ell\right\}  $, so
that $g\geq p+1>p$. Also, $g\in\left\{  p+1,p+2,\ldots,p+\ell\right\}
\subseteq\left\{  1,2,\ldots,p+\ell\right\}  $ (since $p+1\geq1$). Hence,
(\ref{pf.lem.sol.det.pluecker.multi.2.b.gcol}) yields%
\begin{align}
&  \left(  \text{the }g\text{-th column of the matrix }\left(  A\mid
B_{\bullet,I}\right)  \right) \nonumber\\
&  =%
\begin{cases}
A_{\bullet,g}, & \text{if }g\leq p;\\
B_{\bullet,i_{g-p}}, & \text{if }g>p
\end{cases}
=B_{\bullet,i_{g-p}}\ \ \ \ \ \ \ \ \ \ \left(  \text{since }g>p\right)
\nonumber\\
&  =\left(  \text{the }i_{g-p}\text{-th column of the matrix }B\right)
\label{pf.lem.sol.det.pluecker.multi.2.b.cxg2.pf.1}%
\end{align}
(since $B_{\bullet,i_{g-p}}$ is defined to be the $i_{g-p}$-th column of the
matrix $B$). But
\[
\left(  \text{the }\left(  x,g\right)  \text{-th entry of the matrix }\left(
A\mid B_{\bullet,I}\right)  \right)  =c_{x,g}%
\]
(since $\left(  A\mid B_{\bullet,I}\right)  =\left(  c_{i,j}\right)  _{1\leq
i\leq n,\ 1\leq j\leq p+\ell}$). Comparing this with%
\begin{align*}
&  \left(  \text{the }\left(  x,g\right)  \text{-th entry of the matrix
}\left(  A\mid B_{\bullet,I}\right)  \right) \\
&  =\left(  \text{the }x\text{-th entry of }\underbrace{\text{the }g\text{-th
column of the matrix }\left(  A\mid B_{\bullet,I}\right)  }%
_{\substack{=\left(  \text{the }i_{g-p}\text{-th column of the matrix
}B\right)  \\\text{(by (\ref{pf.lem.sol.det.pluecker.multi.2.b.cxg2.pf.1}))}%
}}\right) \\
&  =\left(  \text{the }x\text{-th entry of the }i_{g-p}\text{-th column of the
matrix }B\right) \\
&  =\left(  \text{the }\left(  x,i_{g-p}\right)  \text{-th entry of the matrix
}B\right)  =b_{x,i_{g-p}}\ \ \ \ \ \ \ \ \ \ \left(  \text{since }B=\left(
b_{i,j}\right)  _{1\leq i\leq n,\ 1\leq j\leq m}\right)  ,
\end{align*}
we obtain $c_{x,g}=b_{x,i_{g-p}}$. This proves
(\ref{pf.lem.sol.det.pluecker.multi.2.b.cxg2}).}. Hence, for each
$x\in\left\{  1,2,\ldots,n\right\}  $ and $y\in\left\{  1,2,\ldots
,\ell\right\}  $, we have%
\begin{equation}
c_{x,p+y}=b_{x,i_{y}} \label{pf.lem.sol.det.pluecker.multi.2.b.cxg2b}%
\end{equation}
\footnote{\textit{Proof of (\ref{pf.lem.sol.det.pluecker.multi.2.b.cxg2b}):}
Let $x\in\left\{  1,2,\ldots,n\right\}  $ and $y\in\left\{  1,2,\ldots
,\ell\right\}  $. From $y\in\left\{  1,2,\ldots,\ell\right\}  $, we obtain
$p+y\in\left\{  p+1,p+2,\ldots,p+\ell\right\}  $. Hence,
(\ref{pf.lem.sol.det.pluecker.multi.2.b.cxg2}) (applied to $g=p+y$) yields
$c_{x,p+y}=b_{x,i_{p+y-p}}=b_{x,i_{y}}$ (since $p+y-p=y$). This proves
(\ref{pf.lem.sol.det.pluecker.multi.2.b.cxg2b}).}.

Now, $P$ is a subset of $\left\{  1,2,\ldots,n\right\}  $, whereas $Q$ is a
subset of $\left\{  1,2,\ldots,p+\left\vert I\right\vert \right\}  $ (since
$Q\subseteq\left\{  1,2,\ldots,p+\left\vert I\right\vert \right\}  $). Thus,
the matrix $\operatorname*{sub}\nolimits_{w\left(  P\right)  }^{w\left(
Q\right)  }\left(  A\mid B_{\bullet,I}\right)  $ is well-defined (since
$\left(  A\mid B_{\bullet,I}\right)  $ is an $n\times\left(  p+\left\vert
I\right\vert \right)  $-matrix). Also, $P$ is a subset of $\left\{
1,2,\ldots,n\right\}  $, whereas $Q$ is a subset of $\left\{  1,2,\ldots
,p\right\}  $ (since $Q=\left\{  1,2,\ldots,p\right\}  $). Hence, the matrix
$\operatorname*{sub}\nolimits_{w\left(  P\right)  }^{w\left(  Q\right)  }A$ is
well-defined (since $A$ is an $n\times p$-matrix).

Furthermore, $P$ is a subset of $\left\{  1,2,\ldots,n\right\}  $, whereas $R$
is a subset of $\left\{  1,2,\ldots,p+\left\vert I\right\vert \right\}  $
(since $R\subseteq\left\{  1,2,\ldots,p+\left\vert I\right\vert \right\}  $).
Thus, the matrix $\operatorname*{sub}\nolimits_{w\left(  P\right)  }^{w\left(
R\right)  }\left(  A\mid B_{\bullet,I}\right)  $ is well-defined (since
$\left(  A\mid B_{\bullet,I}\right)  $ is an $n\times\left(  p+\left\vert
I\right\vert \right)  $-matrix). Also, $P$ is a subset of $\left\{
1,2,\ldots,n\right\}  $, whereas $I$ is a subset of $\left\{  1,2,\ldots
,m\right\}  $. Hence, the matrix $\operatorname*{sub}\nolimits_{w\left(
P\right)  }^{w\left(  I\right)  }B$ is well-defined (since $B$ is an $n\times
m$-matrix).

We have $w\left(  Q\right)  =\left(  1,2,\ldots,p\right)  $%
\ \ \ \ \footnote{\textit{Proof.} Recall that $w\left(  Q\right)  $ is the
list of all elements of $Q$ in increasing order (with no repetitions) (by the
definition of $w\left(  Q\right)  $). Thus,%
\begin{align*}
w\left(  Q\right)   &  =\left(  \text{the list of all elements of
}\underbrace{Q}_{=\left\{  1,2,\ldots,p\right\}  }\text{ in increasing order
(with no repetitions)}\right) \\
&  =\left(  \text{the list of all elements of }\left\{  1,2,\ldots,p\right\}
\text{ in increasing order (with no repetitions)}\right) \\
&  =\left(  1,2,\ldots,p\right)  .
\end{align*}
} and $w\left(  R\right)  =\left(  p+1,p+2,\ldots,p+\ell\right)
$\ \ \ \ \footnote{\textit{Proof.} We have $R=\left\{  p+1,p+2,\ldots
,p+\left\vert I\right\vert \right\}  =\left\{  p+1,p+2,\ldots,p+\ell\right\}
$ (since $\left\vert I\right\vert =\ell$).
\par
But $w\left(  R\right)  $ is the list of all elements of $R$ in increasing
order (with no repetitions) (by the definition of $w\left(  R\right)  $).
Thus,%
\begin{align*}
w\left(  R\right)   &  =\left(  \text{the list of all elements of
}\underbrace{R}_{=\left\{  p+1,p+2,\ldots,p+\ell\right\}  }\text{ in
increasing order (with no repetitions)}\right) \\
&  =\left(  \text{the list of all elements of }\left\{  p+1,p+2,\ldots
,p+\ell\right\}  \text{ in increasing order (with no repetitions)}\right) \\
&  =\left(  p+1,p+2,\ldots,p+\ell\right)  .
\end{align*}
}.

Write the list $w\left(  P\right)  $ in the form $w\left(  P\right)  =\left(
p_{1},p_{2},\ldots,p_{h}\right)  $. Proposition
\ref{prop.sect.laplace.notations.w(I)} (applied to $P$ instead of $I$) shows
that $w\left(  P\right)  $ is a $\left\vert P\right\vert $-tuple of elements
of $P$. In other words, $\left(  p_{1},p_{2},\ldots,p_{h}\right)  $ is a
$\left\vert P\right\vert $-tuple of elements of $P$ (since $w\left(  P\right)
=\left(  p_{1},p_{2},\ldots,p_{h}\right)  $). Thus, $p_{1},p_{2},\ldots,p_{h}$
are elements of $P$. Hence, $p_{1},p_{2},\ldots,p_{h}$ are elements of
$\left\{  1,2,\ldots,n\right\}  $ (since $P$ is a subset of $\left\{
1,2,\ldots,n\right\}  $). In other words, $p_{x}\in\left\{  1,2,\ldots
,n\right\}  $ for each $x\in\left\{  1,2,\ldots,h\right\}  $. Therefore, for
each $x\in\left\{  1,2,\ldots,h\right\}  $ and $y\in\left\{  1,2,\ldots
,p\right\}  $, we have%
\begin{equation}
c_{p_{x},y}=a_{p_{x},y} \label{pf.lem.sol.det.pluecker.multi.2.b.cpxy}%
\end{equation}
(by (\ref{pf.lem.sol.det.pluecker.multi.2.b.cxg}) (applied to $p_{x}$ and $y$
instead of $x$ and $g$)).

Recall that $p_{x}\in\left\{  1,2,\ldots,n\right\}  $ for each $x\in\left\{
1,2,\ldots,h\right\}  $. Therefore, for each $x\in\left\{  1,2,\ldots
,h\right\}  $ and $y\in\left\{  1,2,\ldots,\ell\right\}  $, we have%
\begin{equation}
c_{p_{x},p+y}=b_{p_{x},i_{y}} \label{pf.lem.sol.det.pluecker.multi.2.b.cpxy2b}%
\end{equation}
(by (\ref{pf.lem.sol.det.pluecker.multi.2.b.cxg2b}) (applied to $p_{x}$
instead of $x$)).

We have%
\begin{align*}
\operatorname*{sub}\nolimits_{w\left(  P\right)  }^{w\left(  Q\right)
}\left(  A\mid B_{\bullet,I}\right)   &  =\operatorname*{sub}%
\nolimits_{\left(  p_{1},p_{2},\ldots,p_{h}\right)  }^{\left(  1,2,\ldots
,p\right)  }\left(  A\mid B_{\bullet,I}\right)  \ \ \ \ \ \ \ \ \ \ \left(
\begin{array}
[c]{c}%
\text{since }w\left(  P\right)  =\left(  p_{1},p_{2},\ldots,p_{h}\right) \\
\text{and }w\left(  Q\right)  =\left(  1,2,\ldots,p\right)
\end{array}
\right) \\
&  =\operatorname*{sub}\nolimits_{p_{1},p_{2},\ldots,p_{h}}^{1,2,\ldots
,p}\left(  A\mid B_{\bullet,I}\right)  =\left(  \underbrace{c_{p_{x},y}%
}_{\substack{=a_{p_{x},y}\\\text{(by
(\ref{pf.lem.sol.det.pluecker.multi.2.b.cpxy}))}}}\right)  _{1\leq x\leq
h,\ 1\leq y\leq p}\\
&  \ \ \ \ \ \ \ \ \ \ \left(
\begin{array}
[c]{c}%
\text{by the definition of }\operatorname*{sub}\nolimits_{p_{1},p_{2}%
,\ldots,p_{h}}^{1,2,\ldots,p}\left(  A\mid B_{\bullet,I}\right)  \text{,}\\
\text{since }\left(  A\mid B_{\bullet,I}\right)  =\left(  c_{i,j}\right)
_{1\leq i\leq n,\ 1\leq j\leq p+\ell}%
\end{array}
\right) \\
&  =\left(  a_{p_{x},y}\right)  _{1\leq x\leq h,\ 1\leq y\leq p}.
\end{align*}
Comparing this with%
\begin{align*}
\operatorname*{sub}\nolimits_{w\left(  P\right)  }^{w\left(  Q\right)  }A  &
=\operatorname*{sub}\nolimits_{\left(  p_{1},p_{2},\ldots,p_{h}\right)
}^{\left(  1,2,\ldots,p\right)  }A\ \ \ \ \ \ \ \ \ \ \left(
\begin{array}
[c]{c}%
\text{since }w\left(  P\right)  =\left(  p_{1},p_{2},\ldots,p_{h}\right) \\
\text{and }w\left(  Q\right)  =\left(  1,2,\ldots,p\right)
\end{array}
\right) \\
&  =\operatorname*{sub}\nolimits_{p_{1},p_{2},\ldots,p_{h}}^{1,2,\ldots
,p}A=\left(  a_{p_{x},y}\right)  _{1\leq x\leq h,\ 1\leq y\leq p}\\
&  \ \ \ \ \ \ \ \ \ \ \left(
\begin{array}
[c]{c}%
\text{by the definition of }\operatorname*{sub}\nolimits_{p_{1},p_{2}%
,\ldots,p_{h}}^{1,2,\ldots,p}A\text{,}\\
\text{since }A=\left(  a_{i,j}\right)  _{1\leq i\leq n,\ 1\leq j\leq p}%
\end{array}
\right)  ,
\end{align*}
we obtain%
\[
\operatorname*{sub}\nolimits_{w\left(  P\right)  }^{w\left(  Q\right)
}\left(  A\mid B_{\bullet,I}\right)  =\operatorname*{sub}\nolimits_{w\left(
P\right)  }^{w\left(  Q\right)  }A.
\]

Next, we notice that
\begin{align*}
\operatorname*{sub}\nolimits_{w\left(  P\right)  }^{w\left(  R\right)
}\left(  A\mid B_{\bullet,I}\right)   &  =\operatorname*{sub}%
\nolimits_{\left(  p_{1},p_{2},\ldots,p_{h}\right)  }^{\left(  p+1,p+2,\ldots
,p+\ell\right)  }\left(  A\mid B_{\bullet,I}\right) \\
&  \ \ \ \ \ \ \ \ \ \ \left(
\begin{array}
[c]{c}%
\text{since }w\left(  P\right)  =\left(  p_{1},p_{2},\ldots,p_{h}\right) \\
\text{and }w\left(  R\right)  =\left(  p+1,p+2,\ldots,p+\ell\right)
\end{array}
\right) \\
&  =\operatorname*{sub}\nolimits_{p_{1},p_{2},\ldots,p_{h}}^{p+1,p+2,\ldots
,p+\ell}\left(  A\mid B_{\bullet,I}\right)  =\left(  \underbrace{c_{p_{x}%
,p+y}}_{\substack{=b_{p_{x},i_{y}}\\\text{(by
(\ref{pf.lem.sol.det.pluecker.multi.2.b.cpxy2b}))}}}\right)  _{1\leq x\leq
h,\ 1\leq y\leq\ell}\\
&  \ \ \ \ \ \ \ \ \ \ \left(
\begin{array}
[c]{c}%
\text{by the definition of }\operatorname*{sub}\nolimits_{p_{1},p_{2}%
,\ldots,p_{h}}^{p+1,p+2,\ldots,p+\ell}\left(  A\mid B_{\bullet,I}\right)
\text{,}\\
\text{since }\left(  A\mid B_{\bullet,I}\right)  =\left(  c_{i,j}\right)
_{1\leq i\leq n,\ 1\leq j\leq p+\ell}%
\end{array}
\right) \\
&  =\left(  b_{p_{x},i_{y}}\right)  _{1\leq x\leq h,\ 1\leq y\leq\ell}.
\end{align*}
Comparing this with%
\begin{align*}
\operatorname*{sub}\nolimits_{w\left(  P\right)  }^{w\left(  I\right)  }B  &
=\operatorname*{sub}\nolimits_{\left(  p_{1},p_{2},\ldots,p_{h}\right)
}^{\left(  i_{1},i_{2},\ldots,i_{\ell}\right)  }B\ \ \ \ \ \ \ \ \ \ \left(
\begin{array}
[c]{c}%
\text{since }w\left(  I\right)  =\left(  i_{1},i_{2},\ldots,i_{\ell}\right) \\
\text{and }w\left(  P\right)  =\left(  p_{1},p_{2},\ldots,p_{h}\right)
\end{array}
\right) \\
&  =\operatorname*{sub}\nolimits_{p_{1},p_{2},\ldots,p_{h}}^{i_{1}%
,i_{2},\ldots,i_{\ell}}B=\left(  b_{p_{x},i_{y}}\right)  _{1\leq x\leq
h,\ 1\leq y\leq\ell}\\
&  \ \ \ \ \ \ \ \ \ \ \left(
\begin{array}
[c]{c}%
\text{by the definition of }\operatorname*{sub}\nolimits_{p_{1},p_{2}%
,\ldots,p_{h}}^{i_{1},i_{2},\ldots,i_{\ell}}B\text{,}\\
\text{since }B=\left(  b_{i,j}\right)  _{1\leq i\leq n,\ 1\leq j\leq m}%
\end{array}
\right)  ,
\end{align*}
we obtain%
\[
\operatorname*{sub}\nolimits_{w\left(  P\right)  }^{w\left(  R\right)
}\left(  A\mid B_{\bullet,I}\right)  =\operatorname*{sub}\nolimits_{w\left(
P\right)  }^{w\left(  I\right)  }B.
\]
Thus, the proof of Lemma \ref{lem.sol.det.pluecker.multi.2} \textbf{(b)} is complete.
\end{proof}
\end{verlong}

Our next lemma is just a restatement of Exercise
\ref{exe.det.laplace-multi.0r} \textbf{(a)} using more opportune notations:

\begin{lemma}
\label{lem.sol.det.pluecker.multi.3}Let $n\in\mathbb{N}$ and $m\in\mathbb{N}$.
Let $J$ and $K$ be two subsets of $\left\{  1,2,\ldots,m\right\}  $ satisfying
$\left\vert J\right\vert +\left\vert K\right\vert =n$ and $J\cap
K\neq\varnothing$. Let $A\in\mathbb{K}^{m\times n}$. Then,%
\[
\sum_{\substack{I\subseteq\left\{  1,2,\ldots,n\right\}  ;\\\left\vert
I\right\vert =\left\vert J\right\vert }}\left(  -1\right)  ^{\sum I}%
\det\left(  \operatorname*{sub}\nolimits_{w\left(  J\right)  }^{w\left(
I\right)  }A\right)  \det\left(  \operatorname*{sub}\nolimits_{w\left(
K\right)  }^{w\left(  \left[  n\right]  \setminus I\right)  }A\right)  =0.
\]

\end{lemma}

\begin{proof}
[Proof of Lemma \ref{lem.sol.det.pluecker.multi.3}.]We have $\left[  n\right]
=\left\{  1,2,\ldots,n\right\}  $ (by the definition of $\left[  n\right]  $).

For any subset $I$ of $\left\{  1,2,\ldots,n\right\}  $, we let $\widetilde{I}%
$ denote the complement $\left\{  1,2,\ldots,n\right\}  \setminus I$ of $I$.
Thus, for any subset $I$ of $\left\{  1,2,\ldots,n\right\}  $, we have%
\begin{equation}
\widetilde{I}=\underbrace{\left\{  1,2,\ldots,n\right\}  }_{=\left[  n\right]
}\setminus I=\left[  n\right]  \setminus I.
\label{pf.lem.sol.det.pluecker.multi.3.1}%
\end{equation}

Now, Exercise \ref{exe.det.laplace-multi.0r} \textbf{(a)} yields%
\begin{align*}
0  &  =\sum_{\substack{Q\subseteq\left\{  1,2,\ldots,n\right\}  ;\\\left\vert
Q\right\vert =\left\vert J\right\vert }}\left(  -1\right)  ^{\sum Q}%
\det\left(  \operatorname*{sub}\nolimits_{w\left(  J\right)  }^{w\left(
Q\right)  }A\right)  \det\left(  \operatorname*{sub}\nolimits_{w\left(
K\right)  }^{w\left(  \widetilde{Q}\right)  }A\right) \\
&  =\sum_{\substack{I\subseteq\left\{  1,2,\ldots,n\right\}  ;\\\left\vert
I\right\vert =\left\vert J\right\vert }}\left(  -1\right)  ^{\sum I}%
\det\left(  \operatorname*{sub}\nolimits_{w\left(  J\right)  }^{w\left(
I\right)  }A\right)  \underbrace{\det\left(  \operatorname*{sub}%
\nolimits_{w\left(  K\right)  }^{w\left(  \widetilde{I}\right)  }A\right)
}_{\substack{=\det\left(  \operatorname*{sub}\nolimits_{w\left(  K\right)
}^{w\left(  \left[  n\right]  \setminus I\right)  }A\right)  \\\text{(since
}\widetilde{I}=\left[  n\right]  \setminus I\\\text{(by
(\ref{pf.lem.sol.det.pluecker.multi.3.1})))}}}\\
&  \ \ \ \ \ \ \ \ \ \ \left(  \text{here, we have renamed the summation index
}Q\text{ as }I\right) \\
&  =\sum_{\substack{I\subseteq\left\{  1,2,\ldots,n\right\}  ;\\\left\vert
I\right\vert =\left\vert J\right\vert }}\left(  -1\right)  ^{\sum I}%
\det\left(  \operatorname*{sub}\nolimits_{w\left(  J\right)  }^{w\left(
I\right)  }A\right)  \det\left(  \operatorname*{sub}\nolimits_{w\left(
K\right)  }^{w\left(  \left[  n\right]  \setminus I\right)  }A\right)  .
\end{align*}
This proves Lemma \ref{lem.sol.det.pluecker.multi.3}.
\end{proof}

We can now solve Exercise \ref{exe.det.pluecker.multi}:

\begin{proof}
[Solution to Exercise \ref{exe.det.pluecker.multi}.]Combining $k\in\mathbb{N}$
(since $k$ is a positive integer) with $k\leq n$, we obtain $k\in\left\{
0,1,\ldots,n\right\}  $. Hence, $n-k\in\left\{  0,1,\ldots,n\right\}  $.

Let $p=n-k$. Thus, $p=n-k\in\left\{  0,1,\ldots,n\right\}  \subseteq
\mathbb{N}$, so that $p\geq0$ and therefore $p+1\geq0+1=1$. Also, $p\leq n$
(since $p\in\left\{  0,1,\ldots,n\right\}  $).

\begin{vershort}
Let $Q$ be the set $\left\{  1,2,\ldots,p\right\}  $. Thus, $Q\subseteq
\left\{  1,2,\ldots,n\right\}  $ (since $p\leq n$) and $\left\vert
Q\right\vert =p$.
\end{vershort}

\begin{verlong}
Let $Q$ be the set $\left\{  1,2,\ldots,p\right\}  $. Thus, $Q=\left\{
1,2,\ldots,p\right\}  \subseteq\left\{  1,2,\ldots,n\right\}  $ (since $p\leq
n$). In other words, $Q$ is a subset of $\left\{  1,2,\ldots,n\right\}  $.

From $Q=\left\{  1,2,\ldots,p\right\}  $, we obtain $\left\vert Q\right\vert
=\left\vert \left\{  1,2,\ldots,p\right\}  \right\vert =p$.
\end{verlong}

Notice also that $Q=\left\{  1,2,\ldots,p\right\}  =\left\{  1,2,\ldots
,n-k\right\}  $ (since $p=n-k$), so that $Q$ is a subset of $\left\{
1,2,\ldots,n-k\right\}  $. Hence, a matrix $\operatorname*{sub}%
\nolimits_{w\left(  P\right)  }^{w\left(  Q\right)  }A$ is well-defined
whenever $P$ is a subset of $\left\{  1,2,\ldots,n\right\}  $ (since $A$ is an
$n\times\left(  n-k\right)  $-matrix).

\begin{vershort}
Let $R$ be the set $\left\{  p+1,p+2,\ldots,n\right\}  $. Thus, $R\subseteq
\left\{  1,2,\ldots,n\right\}  $ (since $p\geq0$).
\end{vershort}

\begin{verlong}
Let $R$ be the set $\left\{  p+1,p+2,\ldots,n\right\}  $. Thus, $R=\left\{
p+1,p+2,\ldots,n\right\}  \subseteq\left\{  1,2,\ldots,n\right\}  $ (since
$p+1\geq1$). In other words, $R$ is a subset of $\left\{  1,2,\ldots
,n\right\}  $.
\end{verlong}

For any subset $I$ of $\left\{  1,2,\ldots,n\right\}  $, we let $\widetilde{I}%
$ denote the complement $\left\{  1,2,\ldots,n\right\}  \setminus I$ of $I$.
Thus,%
\begin{align*}
\widetilde{Q}  &  =\left\{  1,2,\ldots,n\right\}  \setminus\underbrace{Q}%
_{=\left\{  1,2,\ldots,p\right\}  }=\left\{  1,2,\ldots,n\right\}
\setminus\left\{  1,2,\ldots,p\right\} \\
&  =\left\{  p+1,p+2,\ldots,n\right\}  =R\ \ \ \ \ \ \ \ \ \ \left(
\text{since }R=\left\{  p+1,p+2,\ldots,n\right\}  \right)  .
\end{align*}

For any subset $P$ of $\left\{  1,2,\ldots,n\right\}  $ satisfying $\left\vert
P\right\vert =\left\vert Q\right\vert $, we define an element $\gamma_{P}%
\in\mathbb{K}$ by%
\begin{equation}
\gamma_{P}=\left(  -1\right)  ^{\sum P+\sum Q}\det\left(  \operatorname*{sub}%
\nolimits_{w\left(  P\right)  }^{w\left(  Q\right)  }A\right)  .
\label{sol.det.pluecker.multi.gammaP=}%
\end{equation}

\begin{verlong}
\noindent(It is easy to see that this is
well-defined.\footnote{\textit{Proof.} Let $P$ be a subset of $\left\{
1,2,\ldots,n\right\}  $ satisfying $\left\vert P\right\vert =\left\vert
Q\right\vert $. Then, $P$ is a subset of $\left\{  1,2,\ldots,n\right\}  $,
whereas $Q$ is a subset of $\left\{  1,2,\ldots,n-k\right\}  $ (since
$Q=\left\{  1,2,\ldots,n-k\right\}  $). Thus, the matrix $\operatorname*{sub}%
\nolimits_{w\left(  P\right)  }^{w\left(  Q\right)  }A$ is well-defined (since
$A$ is an $n\times\left(  n-k\right)  $-matrix). This matrix
$\operatorname*{sub}\nolimits_{w\left(  P\right)  }^{w\left(  Q\right)  }A$ is
a $\left\vert P\right\vert \times\left\vert Q\right\vert $-matrix, and thus a
$\left\vert Q\right\vert \times\left\vert Q\right\vert $-matrix (since
$\left\vert P\right\vert =\left\vert Q\right\vert $); hence, it is a square
matrix. Thus, its determinant $\det\left(  \operatorname*{sub}%
\nolimits_{w\left(  P\right)  }^{w\left(  Q\right)  }A\right)  $ is
well-defined as well. This shows that $\gamma_{P}$ is well-defined. Qed.})
\end{verlong}

Next, let us prove that if $I$ is any subset of $\left\{  1,2,\ldots
,n+k\right\}  $ satisfying $\left\vert I\right\vert =k$, then the determinants
$\det\left(  A\mid B_{\bullet,I}\right)  $ and $\det\left(  B_{\bullet,\sim
I}\right)  $ are well-defined and the equality
\begin{equation}
\det\left(  A\mid B_{\bullet,I}\right)  =\sum_{\substack{P\subseteq\left\{
1,2,\ldots,n\right\}  ;\\\left\vert P\right\vert =\left\vert Q\right\vert
}}\gamma_{P}\det\left(  \operatorname*{sub}\nolimits_{w\left(  \widetilde{P}%
\right)  }^{w\left(  I\right)  }B\right)  \label{sol.det.pluecker.multi.step2}%
\end{equation}
holds.

[\textit{Proof:} Let $I$ be a subset of $\left\{  1,2,\ldots,n+k\right\}  $
satisfying $\left\vert I\right\vert =k$. Then, $\left(  A\mid B_{\bullet
,I}\right)  $ is an $n\times\left(  n-k+\left\vert I\right\vert \right)
$-matrix (by the definition of $\left(  A\mid B_{\bullet,I}\right)  $). In
other words, $\left(  A\mid B_{\bullet,I}\right)  $ is an $n\times n$-matrix
(since $n-k+\underbrace{\left\vert I\right\vert }_{=k}=n-k+k=n$). Thus, its
determinant $\det\left(  A\mid B_{\bullet,I}\right)  $ is well-defined.

Also, $B_{\bullet,\sim I}$ is an $n\times\left(  n+k-\left\vert I\right\vert
\right)  $-matrix (by the definition of $B_{\bullet,\sim I}$). In other words,
$B_{\bullet,\sim I}$ is an $n\times n$-matrix (since
$n+k-\underbrace{\left\vert I\right\vert }_{=k}=n+k-k=n$). Hence, its
determinant $\det\left(  B_{\bullet,\sim I}\right)  $ is well-defined.

It remains to show that the equality (\ref{sol.det.pluecker.multi.step2}) holds.

We have $A\in\mathbb{K}^{n\times\left(  n-k\right)  }=\mathbb{K}^{n\times p}$
(since $n-k=p$). Furthermore, $R$ is the set $\left\{  p+1,p+2,\ldots
,p+\left\vert I\right\vert \right\}  $\ \ \ \ \footnote{\textit{Proof.} Adding
the equalities $p=n-k$ and $\left\vert I\right\vert =k$ together, we obtain
$p+\left\vert I\right\vert =\left(  n-k\right)  +k=n$. Hence, $\left\{
p+1,p+2,\ldots,p+\left\vert I\right\vert \right\}  =\left\{  p+1,p+2,\ldots
,n\right\}  $. Comparing this with $R=\left\{  p+1,p+2,\ldots,n\right\}  $, we
obtain $R=\left\{  p+1,p+2,\ldots,p+\left\vert I\right\vert \right\}  $. In
other words, $R$ is the set $\left\{  p+1,p+2,\ldots,p+\left\vert I\right\vert
\right\}  $.}. Hence, Lemma \ref{lem.sol.det.pluecker.multi.2} \textbf{(b)}
(applied to $m=n+k$) shows that if $P$ is any subset of $\left\{
1,2,\ldots,n\right\}  $, then%
\begin{equation}
\operatorname*{sub}\nolimits_{w\left(  P\right)  }^{w\left(  Q\right)
}\left(  A\mid B_{\bullet,I}\right)  =\operatorname*{sub}\nolimits_{w\left(
P\right)  }^{w\left(  Q\right)  }A \label{sol.det.pluecker.multi.l2u1}%
\end{equation}
and%
\begin{equation}
\operatorname*{sub}\nolimits_{w\left(  P\right)  }^{w\left(  R\right)
}\left(  A\mid B_{\bullet,I}\right)  =\operatorname*{sub}\nolimits_{w\left(
P\right)  }^{w\left(  I\right)  }B. \label{sol.det.pluecker.multi.l2u2}%
\end{equation}

Also, if $P$ is any subset of $\left\{  1,2,\ldots,n\right\}  $, then
$\widetilde{P}$ is a subset of $\left\{  1,2,\ldots,n\right\}  $ as well
(since the definition of $\widetilde{P}$ yields $\widetilde{P}=\left\{
1,2,\ldots,n\right\}  \setminus P\subseteq\left\{  1,2,\ldots,n\right\}  $).

But $\left(  A\mid B_{\bullet,I}\right)  $ is an $n\times n$-matrix. In other
words, $\left(  A\mid B_{\bullet,I}\right)  \in\mathbb{K}^{n\times n}$. Hence,
Theorem \ref{thm.det.laplace-multi} \textbf{(b)} (applied to $\left(  A\mid
B_{\bullet,I}\right)  $ instead of $A$) yields%
\begin{align*}
&  \det\left(  A\mid B_{\bullet,I}\right) \\
&  =\sum_{\substack{P\subseteq\left\{  1,2,\ldots,n\right\}  ;\\\left\vert
P\right\vert =\left\vert Q\right\vert }}\left(  -1\right)  ^{\sum P+\sum
Q}\det\left(  \underbrace{\operatorname*{sub}\nolimits_{w\left(  P\right)
}^{w\left(  Q\right)  }\left(  A\mid B_{\bullet,I}\right)  }%
_{\substack{=\operatorname*{sub}\nolimits_{w\left(  P\right)  }^{w\left(
Q\right)  }A\\\text{(by (\ref{sol.det.pluecker.multi.l2u1}))}}}\right)
\underbrace{\det\left(  \operatorname*{sub}\nolimits_{w\left(  \widetilde{P}%
\right)  }^{w\left(  \widetilde{Q}\right)  }\left(  A\mid B_{\bullet
,I}\right)  \right)  }_{\substack{=\det\left(  \operatorname*{sub}%
\nolimits_{w\left(  \widetilde{P}\right)  }^{w\left(  R\right)  }\left(  A\mid
B_{\bullet,I}\right)  \right)  \\\text{(since }\widetilde{Q}=R\text{)}}}\\
&  =\sum_{\substack{P\subseteq\left\{  1,2,\ldots,n\right\}  ;\\\left\vert
P\right\vert =\left\vert Q\right\vert }}\underbrace{\left(  -1\right)  ^{\sum
P+\sum Q}\det\left(  \operatorname*{sub}\nolimits_{w\left(  P\right)
}^{w\left(  Q\right)  }A\right)  }_{\substack{=\gamma_{P}\\\text{(by
(\ref{sol.det.pluecker.multi.gammaP=}))}}}\det\left(
\underbrace{\operatorname*{sub}\nolimits_{w\left(  \widetilde{P}\right)
}^{w\left(  R\right)  }\left(  A\mid B_{\bullet,I}\right)  }%
_{\substack{=\operatorname*{sub}\nolimits_{w\left(  \widetilde{P}\right)
}^{w\left(  I\right)  }B\\\text{(by (\ref{sol.det.pluecker.multi.l2u2}%
)}\\\text{(applied to }\widetilde{P}\text{ instead of }P\text{))}}}\right) \\
&  =\sum_{\substack{P\subseteq\left\{  1,2,\ldots,n\right\}  ;\\\left\vert
P\right\vert =\left\vert Q\right\vert }}\gamma_{P}\det\left(
\operatorname*{sub}\nolimits_{w\left(  \widetilde{P}\right)  }^{w\left(
I\right)  }B\right)  .
\end{align*}
In other words, the equality (\ref{sol.det.pluecker.multi.step2}) holds.

Now, forget that we fixed $I$. We thus have shown that if $I$ is any subset of
$\left\{  1,2,\ldots,n+k\right\}  $ satisfying $\left\vert I\right\vert =k$,
then the determinants $\det\left(  A\mid B_{\bullet,I}\right)  $ and
$\det\left(  B_{\bullet,\sim I}\right)  $ are well-defined and the equality
(\ref{sol.det.pluecker.multi.step2}) holds. Qed.]

On the other hand, if $P$ is any subset of $\left\{  1,2,\ldots,n\right\}  $
satisfying $\left\vert P\right\vert =\left\vert Q\right\vert $, then%
\begin{equation}
\sum_{\substack{I\subseteq\left\{  1,2,\ldots,n+k\right\}  ;\\\left\vert
I\right\vert =k}}\left(  -1\right)  ^{\sum I}\det\left(  \operatorname*{sub}%
\nolimits_{w\left(  \widetilde{P}\right)  }^{w\left(  I\right)  }B\right)
\det\left(  B_{\bullet,\sim I}\right)  =0.
\label{sol.det.pluecker.multi.step3}%
\end{equation}

[\textit{Proof of (\ref{sol.det.pluecker.multi.step3}):} Let $P$ be any subset
of $\left\{  1,2,\ldots,n\right\}  $ satisfying $\left\vert P\right\vert
=\left\vert Q\right\vert $. Thus, $\left\vert P\right\vert =\left\vert
Q\right\vert =p=n-k$.

\begin{vershort}
Since $P$ is a subset of $\left\{  1,2,\ldots,n\right\}  $, its complement
$\widetilde{P}$ is a subset of $\left\{  1,2,\ldots,n\right\}  $ as well and
satisfies%
\[
\left\vert \widetilde{P}\right\vert =\underbrace{\left\vert \left\{
1,2,\ldots,n\right\}  \right\vert }_{=n}-\underbrace{\left\vert P\right\vert
}_{=n-k}=n-\left(  n-k\right)  =k.
\]

\end{vershort}

\begin{verlong}
The definition of $\widetilde{P}$ yields $\widetilde{P}=\left\{
1,2,\ldots,n\right\}  \setminus P\subseteq\left\{  1,2,\ldots,n\right\}  $, so
that $\widetilde{P}$ is a subset of $\left\{  1,2,\ldots,n\right\}  $. Also,
from $\widetilde{P}=\left\{  1,2,\ldots,n\right\}  \setminus P$, we obtain%
\begin{align*}
\left\vert \widetilde{P}\right\vert  &  =\left\vert \left\{  1,2,\ldots
,n\right\}  \setminus P\right\vert =\underbrace{\left\vert \left\{
1,2,\ldots,n\right\}  \right\vert }_{=n}-\underbrace{\left\vert P\right\vert
}_{=n-k}\ \ \ \ \ \ \ \ \ \ \left(  \text{since }P\text{ is a subset of
}\left\{  1,2,\ldots,n\right\}  \right) \\
&  =n-\left(  n-k\right)  =k.
\end{align*}

\end{verlong}

\begin{vershort}
Furthermore, $\left[  n\right]  $ is a subset of $\left\{  1,2,\ldots
,n\right\}  $ (since $\left[  n\right]  =\left\{  1,2,\ldots,n\right\}  $) and
satisfies $\underbrace{\left\vert \widetilde{P}\right\vert }_{=k}%
+\underbrace{\left\vert \left[  n\right]  \right\vert }_{=n}=k+n=n+k$. Also,
$\left\vert \widetilde{P}\right\vert =k>0$ (since $k$ is a positive integer),
so that $\widetilde{P}\neq\varnothing$. But $\widetilde{P}\subseteq\left\{
1,2,\ldots,n\right\}  =\left[  n\right]  $ and thus $\widetilde{P}\cap\left[
n\right]  =\widetilde{P}\neq\varnothing$.
\end{vershort}

\begin{verlong}
Furthermore, $\left[  n\right]  $ is a subset of $\left\{  1,2,\ldots
,n\right\}  $ (since $\left[  n\right]  =\left\{  1,2,\ldots,n\right\}  $).
Clearly,
\[
\underbrace{\left\vert \widetilde{P}\right\vert }_{=k}+\left\vert
\underbrace{\left[  n\right]  }_{=\left\{  1,2,\ldots,n\right\}  }\right\vert
=k+\underbrace{\left\vert \left\{  1,2,\ldots,n\right\}  \right\vert }%
_{=n}=k+n=n+k.
\]
Also, $\left\vert \widetilde{P}\right\vert =k>0$ (since $k$ is a positive
integer), so that $\widetilde{P}\neq\varnothing$. But $\widetilde{P}%
\subseteq\left\{  1,2,\ldots,n\right\}  =\left[  n\right]  $ and thus
$\widetilde{P}\cap\left[  n\right]  =\widetilde{P}\neq\varnothing$.
\end{verlong}

Hence, Lemma \ref{lem.sol.det.pluecker.multi.3} (applied to $n$, $n+k$, $B$,
$\widetilde{P}$ and $\left[  n\right]  $ instead of $m$, $n$, $A$, $J$ and
$K$) yields%
\[
\sum_{\substack{I\subseteq\left\{  1,2,\ldots,n+k\right\}  ;\\\left\vert
I\right\vert =\left\vert \widetilde{P}\right\vert }}\left(  -1\right)  ^{\sum
I}\det\left(  \operatorname*{sub}\nolimits_{w\left(  \widetilde{P}\right)
}^{w\left(  I\right)  }B\right)  \det\left(  \operatorname*{sub}%
\nolimits_{w\left(  \left[  n\right]  \right)  }^{w\left(  \left[  n+k\right]
\setminus I\right)  }B\right)  =0.
\]
In view of $\left\vert \widetilde{P}\right\vert =k$, this rewrites as%
\begin{equation}
\sum_{\substack{I\subseteq\left\{  1,2,\ldots,n+k\right\}  ;\\\left\vert
I\right\vert =k}}\left(  -1\right)  ^{\sum I}\det\left(  \operatorname*{sub}%
\nolimits_{w\left(  \widetilde{P}\right)  }^{w\left(  I\right)  }B\right)
\det\left(  \operatorname*{sub}\nolimits_{w\left(  \left[  n\right]  \right)
}^{w\left(  \left[  n+k\right]  \setminus I\right)  }B\right)  =0.
\label{sol.det.pluecker.multi.step3.pf.3}%
\end{equation}

Now,%
\begin{align*}
&  \sum_{\substack{I\subseteq\left\{  1,2,\ldots,n+k\right\}  ;\\\left\vert
I\right\vert =k}}\left(  -1\right)  ^{\sum I}\det\left(  \operatorname*{sub}%
\nolimits_{w\left(  \widetilde{P}\right)  }^{w\left(  I\right)  }B\right)
\det\left(  \underbrace{B_{\bullet,\sim I}}_{\substack{=\operatorname*{sub}%
\nolimits_{w\left(  \left[  n\right]  \right)  }^{w\left(  \left[  n+k\right]
\setminus I\right)  }B\\\text{(by Lemma \ref{lem.sol.det.pluecker.multi.1}
\textbf{(b)} (applied to }m=n+k\text{))}}}\right) \\
&  =\sum_{\substack{I\subseteq\left\{  1,2,\ldots,n+k\right\}  ;\\\left\vert
I\right\vert =k}}\left(  -1\right)  ^{\sum I}\det\left(  \operatorname*{sub}%
\nolimits_{w\left(  \widetilde{P}\right)  }^{w\left(  I\right)  }B\right)
\det\left(  \operatorname*{sub}\nolimits_{w\left(  \left[  n\right]  \right)
}^{w\left(  \left[  n+k\right]  \setminus I\right)  }B\right)  =0
\end{align*}
(by (\ref{sol.det.pluecker.multi.step3.pf.3})). This proves
(\ref{sol.det.pluecker.multi.step3}).]

\begin{vershort}
Now,%
\begin{align}
&  \sum_{\substack{I\subseteq\left\{  1,2,\ldots,n+k\right\}  ;\\\left\vert
I\right\vert =k}}\left(  -1\right)  ^{\sum I}\underbrace{\det\left(  A\mid
B_{\bullet,I}\right)  }_{\substack{=\sum_{\substack{P\subseteq\left\{
1,2,\ldots,n\right\}  ;\\\left\vert P\right\vert =\left\vert Q\right\vert
}}\gamma_{P}\det\left(  \operatorname*{sub}\nolimits_{w\left(  \widetilde{P}%
\right)  }^{w\left(  I\right)  }B\right)  \\\text{(by
(\ref{sol.det.pluecker.multi.step2}))}}}\det\left(  B_{\bullet,\sim I}\right)
\nonumber\\
&  =\sum_{\substack{I\subseteq\left\{  1,2,\ldots,n+k\right\}  ;\\\left\vert
I\right\vert =k}}\left(  -1\right)  ^{\sum I}\left(  \sum
_{\substack{P\subseteq\left\{  1,2,\ldots,n\right\}  ;\\\left\vert
P\right\vert =\left\vert Q\right\vert }}\gamma_{P}\det\left(
\operatorname*{sub}\nolimits_{w\left(  \widetilde{P}\right)  }^{w\left(
I\right)  }B\right)  \right)  \det\left(  B_{\bullet,\sim I}\right)
\nonumber\\
&  =\sum_{\substack{P\subseteq\left\{  1,2,\ldots,n\right\}  ;\\\left\vert
P\right\vert =\left\vert Q\right\vert }}\gamma_{P}\underbrace{\sum
_{\substack{I\subseteq\left\{  1,2,\ldots,n+k\right\}  ;\\\left\vert
I\right\vert =k}}\left(  -1\right)  ^{\sum I}\det\left(  \operatorname*{sub}%
\nolimits_{w\left(  \widetilde{P}\right)  }^{w\left(  I\right)  }B\right)
\det\left(  B_{\bullet,\sim I}\right)  }_{\substack{=0\\\text{(by
(\ref{sol.det.pluecker.multi.step3}))}}}\nonumber\\
&  =\sum_{\substack{P\subseteq\left\{  1,2,\ldots,n\right\}  ;\\\left\vert
P\right\vert =\left\vert Q\right\vert }}\gamma_{P}0=0.
\label{sol.det.pluecker.multi.short.step5}%
\end{align}
Hence,%
\begin{align*}
&  \sum_{\substack{I\subseteq\left\{  1,2,\ldots,n+k\right\}  ;\\\left\vert
I\right\vert =k}}\underbrace{\left(  -1\right)  ^{\sum I+\left(
1+2+\cdots+k\right)  }}_{=\left(  -1\right)  ^{\sum I}\left(  -1\right)
^{1+2+\cdots+k}}\det\left(  A\mid B_{\bullet,I}\right)  \det\left(
B_{\bullet,\sim I}\right) \\
&  =\left(  -1\right)  ^{1+2+\cdots+k}\underbrace{\sum_{\substack{I\subseteq
\left\{  1,2,\ldots,n+k\right\}  ;\\\left\vert I\right\vert =k}}\left(
-1\right)  ^{\sum I}\det\left(  A\mid B_{\bullet,I}\right)  \det\left(
B_{\bullet,\sim I}\right)  }_{\substack{=0\\\text{(by
(\ref{sol.det.pluecker.multi.short.step5}))}}}\\
&  =0.
\end{align*}
This solves Exercise \ref{exe.det.pluecker.multi}. \qedhere

\end{vershort}

\begin{verlong}
Now,%
\begin{align}
&  \sum_{\substack{I\subseteq\left\{  1,2,\ldots,n+k\right\}  ;\\\left\vert
I\right\vert =k}}\left(  -1\right)  ^{\sum I}\underbrace{\det\left(  A\mid
B_{\bullet,I}\right)  }_{\substack{=\sum_{\substack{P\subseteq\left\{
1,2,\ldots,n\right\}  ;\\\left\vert P\right\vert =\left\vert Q\right\vert
}}\gamma_{P}\det\left(  \operatorname*{sub}\nolimits_{w\left(  \widetilde{P}%
\right)  }^{w\left(  I\right)  }B\right)  \\\text{(by
(\ref{sol.det.pluecker.multi.step2}))}}}\det\left(  B_{\bullet,\sim I}\right)
\nonumber\\
&  =\sum_{\substack{I\subseteq\left\{  1,2,\ldots,n+k\right\}  ;\\\left\vert
I\right\vert =k}}\left(  -1\right)  ^{\sum I}\left(  \sum
_{\substack{P\subseteq\left\{  1,2,\ldots,n\right\}  ;\\\left\vert
P\right\vert =\left\vert Q\right\vert }}\gamma_{P}\det\left(
\operatorname*{sub}\nolimits_{w\left(  \widetilde{P}\right)  }^{w\left(
I\right)  }B\right)  \right)  \det\left(  B_{\bullet,\sim I}\right)
\nonumber\\
&  =\sum_{\substack{I\subseteq\left\{  1,2,\ldots,n+k\right\}  ;\\\left\vert
I\right\vert =k}}\left(  -1\right)  ^{\sum I}\sum_{\substack{P\subseteq
\left\{  1,2,\ldots,n\right\}  ;\\\left\vert P\right\vert =\left\vert
Q\right\vert }}\gamma_{P}\det\left(  \operatorname*{sub}\nolimits_{w\left(
\widetilde{P}\right)  }^{w\left(  I\right)  }B\right)  \det\left(
B_{\bullet,\sim I}\right) \nonumber\\
&  =\underbrace{\sum_{\substack{I\subseteq\left\{  1,2,\ldots,n+k\right\}
;\\\left\vert I\right\vert =k}}\sum_{\substack{P\subseteq\left\{
1,2,\ldots,n\right\}  ;\\\left\vert P\right\vert =\left\vert Q\right\vert }%
}}_{=\sum_{\substack{P\subseteq\left\{  1,2,\ldots,n\right\}  ;\\\left\vert
P\right\vert =\left\vert Q\right\vert }}\sum_{\substack{I\subseteq\left\{
1,2,\ldots,n+k\right\}  ;\\\left\vert I\right\vert =k}}}\left(  -1\right)
^{\sum I}\gamma_{P}\det\left(  \operatorname*{sub}\nolimits_{w\left(
\widetilde{P}\right)  }^{w\left(  I\right)  }B\right)  \det\left(
B_{\bullet,\sim I}\right) \nonumber\\
&  =\sum_{\substack{P\subseteq\left\{  1,2,\ldots,n\right\}  ;\\\left\vert
P\right\vert =\left\vert Q\right\vert }}\underbrace{\sum_{\substack{I\subseteq
\left\{  1,2,\ldots,n+k\right\}  ;\\\left\vert I\right\vert =k}}\left(
-1\right)  ^{\sum I}\gamma_{P}\det\left(  \operatorname*{sub}%
\nolimits_{w\left(  \widetilde{P}\right)  }^{w\left(  I\right)  }B\right)
\det\left(  B_{\bullet,\sim I}\right)  }_{=\gamma_{P}\sum
_{\substack{I\subseteq\left\{  1,2,\ldots,n+k\right\}  ;\\\left\vert
I\right\vert =k}}\left(  -1\right)  ^{\sum I}\det\left(  \operatorname*{sub}%
\nolimits_{w\left(  \widetilde{P}\right)  }^{w\left(  I\right)  }B\right)
\det\left(  B_{\bullet,\sim I}\right)  }\nonumber\\
&  =\sum_{\substack{P\subseteq\left\{  1,2,\ldots,n\right\}  ;\\\left\vert
P\right\vert =\left\vert Q\right\vert }}\gamma_{P}\underbrace{\sum
_{\substack{I\subseteq\left\{  1,2,\ldots,n+k\right\}  ;\\\left\vert
I\right\vert =k}}\left(  -1\right)  ^{\sum I}\det\left(  \operatorname*{sub}%
\nolimits_{w\left(  \widetilde{P}\right)  }^{w\left(  I\right)  }B\right)
\det\left(  B_{\bullet,\sim I}\right)  }_{\substack{=0\\\text{(by
(\ref{sol.det.pluecker.multi.step3}))}}}\nonumber\\
&  =\sum_{\substack{P\subseteq\left\{  1,2,\ldots,n\right\}  ;\\\left\vert
P\right\vert =\left\vert Q\right\vert }}\gamma_{P}0=0.
\label{sol.det.pluecker.multi.step5}%
\end{align}
Hence,%
\begin{align*}
&  \sum_{\substack{I\subseteq\left\{  1,2,\ldots,n+k\right\}  ;\\\left\vert
I\right\vert =k}}\underbrace{\left(  -1\right)  ^{\sum I+\left(
1+2+\cdots+k\right)  }}_{=\left(  -1\right)  ^{\sum I}\left(  -1\right)
^{1+2+\cdots+k}}\det\left(  A\mid B_{\bullet,I}\right)  \det\left(
B_{\bullet,\sim I}\right) \\
&  =\sum_{\substack{I\subseteq\left\{  1,2,\ldots,n+k\right\}  ;\\\left\vert
I\right\vert =k}}\left(  -1\right)  ^{\sum I}\left(  -1\right)  ^{1+2+\cdots
+k}\det\left(  A\mid B_{\bullet,I}\right)  \det\left(  B_{\bullet,\sim
I}\right) \\
&  =\left(  -1\right)  ^{1+2+\cdots+k}\underbrace{\sum_{\substack{I\subseteq
\left\{  1,2,\ldots,n+k\right\}  ;\\\left\vert I\right\vert =k}}\left(
-1\right)  ^{\sum I}\det\left(  A\mid B_{\bullet,I}\right)  \det\left(
B_{\bullet,\sim I}\right)  }_{\substack{=0\\\text{(by
(\ref{sol.det.pluecker.multi.step5}))}}}\\
&  =0.
\end{align*}
This solves Exercise \ref{exe.det.pluecker.multi}.
\end{verlong}
\end{proof}

\subsection{Solution to Exercise \ref{exe.det.rk1upd}}

In this section, we shall use the following notation:

\begin{definition}
\label{def.sol.det.rk1upd.scalar}If $B$ is any $1\times1$-matrix, then
$\operatorname*{ent}B$ will denote the $\left(  1,1\right)  $-th entry of $B$.
(This entry is, of course, the only entry of $B$. Thus, the $1\times1$-matrix
$B$ satisfies $B=\left(
\begin{array}
[c]{c}%
\operatorname*{ent}B
\end{array}
\right)  $.)
\end{definition}

We can now restate Exercise \ref{exe.det.rk1upd} in a form that uses no abuse
of notation (such as identifying $1\times1$-matrices with elements of
$\mathbb{K}$):

\begin{theorem}
\label{thm.sol.det.rk1upd.claim}Let $n\in\mathbb{N}$. Let $u$ be a column
vector with $n$ entries, and let $v$ be a row vector with $n$ entries. (Thus,
$uv$ is an $n\times n$-matrix, whereas $vu$ is a $1\times1$-matrix.) Let $A$
be an $n\times n$-matrix. Then,%
\[
\det\left(  A+uv\right)  =\det A+\operatorname*{ent}\left(  v\left(
\operatorname*{adj}A\right)  u\right)  .
\]

\end{theorem}

Before we prove this theorem, let us make some preparations. First comes a
simple fact:

\begin{proposition}
\label{prop.sol.det.rk1upd.claim.ent(vBu)}Let $n\in\mathbb{N}$ and
$m\in\mathbb{N}$. Let $u=\left(  u_{1},u_{2},\ldots,u_{m}\right)  ^{T}%
\in\mathbb{K}^{m\times1}$ and $v=\left(  v_{1},v_{2},\ldots,v_{n}\right)
\in\mathbb{K}^{1\times n}$. Let $B=\left(  b_{i,j}\right)  _{1\leq i\leq
n,\ 1\leq j\leq m}\in\mathbb{K}^{n\times m}$. Then,%
\[
\operatorname*{ent}\left(  vBu\right)  =\sum_{i=1}^{n}\sum_{j=1}^{m}u_{j}%
v_{i}b_{i,j}.
\]

\end{proposition}

\begin{proof}
[Proof of Proposition \ref{prop.sol.det.rk1upd.claim.ent(vBu)}.]We have
\[
u=\left(  u_{1},u_{2},\ldots,u_{m}\right)  ^{T}=\left(
\begin{array}
[c]{c}%
u_{1}\\
u_{2}\\
\vdots\\
u_{m}%
\end{array}
\right)  =\left(  u_{i}\right)  _{1\leq i\leq m,\ 1\leq j\leq1}%
\]
and%
\[
v=\left(  v_{1},v_{2},\ldots,v_{n}\right)  =\left(  v_{j}\right)  _{1\leq
i\leq1,\ 1\leq j\leq n}.
\]
The definition of the product of two matrices yields%
\[
vB=\left(  \sum_{k=1}^{n}v_{k}b_{k,j}\right)  _{1\leq i\leq1,\ 1\leq j\leq m}%
\]
(since $v=\left(  v_{j}\right)  _{1\leq i\leq1,\ 1\leq j\leq n}$ and
$B=\left(  b_{i,j}\right)  _{1\leq i\leq n,\ 1\leq j\leq m}$). Thus,%
\[
vB=\left(  \underbrace{\sum_{k=1}^{n}v_{k}b_{k,j}}_{\substack{=\sum_{p=1}%
^{n}v_{p}b_{p,j}\\\text{(here, we have renamed the}\\\text{summation index
}k\text{ as }p\text{)}}}\right)  _{1\leq i\leq1,\ 1\leq j\leq m}=\left(
\sum_{p=1}^{n}v_{p}b_{p,j}\right)  _{1\leq i\leq1,\ 1\leq j\leq m}.
\]

Now, the definition of the product of two matrices yields%
\begin{equation}
\left(  vB\right)  u=\left(  \sum_{k=1}^{m}\left(  \sum_{p=1}^{n}v_{p}%
b_{p,k}\right)  u_{k}\right)  _{1\leq i\leq1,\ 1\leq j\leq1}
\label{pf.prop.sol.det.rk1upd.claim.ent(vBu).1}%
\end{equation}
(since $vB=\left(  \sum_{p=1}^{n}v_{p}b_{p,j}\right)  _{1\leq i\leq1,\ 1\leq
j\leq m}$ and $u=\left(  u_{i}\right)  _{1\leq i\leq m,\ 1\leq j\leq1}$). In
particular, $\left(  vB\right)  u$ is a $1\times1$-matrix. The definition of
$\operatorname*{ent}\left(  \left(  vB\right)  u\right)  $ shows that
$\operatorname*{ent}\left(  \left(  vB\right)  u\right)  $ is the $\left(
1,1\right)  $-th entry of this matrix $\left(  vB\right)  u$. Thus,%
\begin{align*}
\operatorname*{ent}\left(  \left(  vB\right)  u\right)   &  =\left(  \text{the
}\left(  1,1\right)  \text{-th entry of the matrix }\underbrace{\left(
vB\right)  u}_{\substack{=\left(  \sum_{k=1}^{m}\left(  \sum_{p=1}^{n}%
v_{p}b_{p,k}\right)  u_{k}\right)  _{1\leq i\leq1,\ 1\leq j\leq1}\\\text{(by
(\ref{pf.prop.sol.det.rk1upd.claim.ent(vBu).1}))}}}\right) \\
&  =\left(  \text{the }\left(  1,1\right)  \text{-th entry of the matrix
}\left(  \sum_{k=1}^{m}\left(  \sum_{p=1}^{n}v_{p}b_{p,k}\right)
u_{k}\right)  _{1\leq i\leq1,\ 1\leq j\leq1}\right) \\
&  =\sum_{k=1}^{m}\left(  \sum_{p=1}^{n}v_{p}b_{p,k}\right)  u_{k}%
=\underbrace{\sum_{k=1}^{m}\sum_{p=1}^{n}}_{=\sum_{p=1}^{n}\sum_{k=1}^{m}%
}\underbrace{v_{p}b_{p,k}u_{k}}_{=u_{k}v_{p}b_{p,k}}\\
&  =\sum_{p=1}^{n}\sum_{k=1}^{m}u_{k}v_{p}b_{p,k}=\sum_{i=1}^{n}\sum_{k=1}%
^{m}u_{k}v_{i}b_{i,k}\\
&  \ \ \ \ \ \ \ \ \ \ \left(  \text{here, we have renamed the summation index
}p\text{ as }i\right) \\
&  =\sum_{i=1}^{n}\sum_{j=1}^{m}u_{j}v_{i}b_{i,j}%
\end{align*}
(here, we have renamed the summation index $k$ as $j$). Since $\left(
vB\right)  u=vBu$, this rewrites as%
\[
\operatorname*{ent}\left(  vBu\right)  =\sum_{i=1}^{n}\sum_{j=1}^{m}u_{j}%
v_{i}b_{i,j}.
\]
This proves Proposition \ref{prop.sol.det.rk1upd.claim.ent(vBu)}.
\end{proof}

From Proposition \ref{prop.sol.det.rk1upd.claim.ent(vBu)}, we can easily
obtain the following:

\begin{corollary}
\label{cor.sol.det.rk1upd.claim.ent(v(adjA)u)}Let $n\in\mathbb{N}$. Let
$u=\left(  u_{1},u_{2},\ldots,u_{n}\right)  ^{T}\in\mathbb{K}^{n\times1}$ and
$v=\left(  v_{1},v_{2},\ldots,v_{n}\right)  \in\mathbb{K}^{1\times n}$. Let
$A$ be an $n\times n$-matrix. Then,%
\[
\operatorname*{ent}\left(  v\left(  \operatorname*{adj}A\right)  u\right)
=\sum_{i=1}^{n}\sum_{j=1}^{n}\left(  -1\right)  ^{i+j}u_{j}v_{i}\det\left(
A_{\sim j,\sim i}\right)  .
\]
(Here, we are using the notations introduced in Definition
\ref{def.submatrix.minor}.)
\end{corollary}

\begin{proof}
[Proof of Corollary \ref{cor.sol.det.rk1upd.claim.ent(v(adjA)u)}.]The
definition of $\operatorname*{adj}A$ yields%
\[
\operatorname*{adj}A=\left(  \left(  -1\right)  ^{i+j}\det\left(  A_{\sim
j,\sim i}\right)  \right)  _{1\leq i\leq n,\ 1\leq j\leq n}.
\]
Hence, Proposition \ref{prop.sol.det.rk1upd.claim.ent(vBu)} (applied to $m=n$,
$B=\operatorname*{adj}A$ and $b_{i,j}=\left(  -1\right)  ^{i+j}\det\left(
A_{\sim j,\sim i}\right)  $) yields%
\begin{align*}
\operatorname*{ent}\left(  v\left(  \operatorname*{adj}A\right)  u\right)   &
=\sum_{i=1}^{n}\sum_{j=1}^{n}\underbrace{u_{j}v_{i}\left(  -1\right)  ^{i+j}%
}_{=\left(  -1\right)  ^{i+j}u_{j}v_{i}}\det\left(  A_{\sim j,\sim i}\right)
\\
&  =\sum_{i=1}^{n}\sum_{j=1}^{n}\left(  -1\right)  ^{i+j}u_{j}v_{i}\det\left(
A_{\sim j,\sim i}\right)  .
\end{align*}
This proves Corollary \ref{cor.sol.det.rk1upd.claim.ent(v(adjA)u)}.
\end{proof}

Next, we state something slightly more interesting:

\begin{proposition}
\label{prop.sol.det.rk1upd.claim.AB+uv}Let $n\in\mathbb{N}$ and $m\in
\mathbb{N}$. Let $A\in\mathbb{K}^{n\times m}$ and $B\in\mathbb{K}^{n\times m}%
$. Let $u\in\mathbb{K}^{n\times1}$ and $v\in\mathbb{K}^{n\times1}$. Then,%
\[
\left(  A\mid u\right)  \left(  B\mid v\right)  ^{T}=AB^{T}+uv^{T}.
\]
(Here, we are using the notations introduced in Definition \ref{def.addcol}.)
\end{proposition}

\begin{example}
If we set $n=2$, $m=3$, $A=\left(
\begin{array}
[c]{ccc}%
a & a^{\prime} & a^{\prime\prime}\\
b & b^{\prime} & b^{\prime\prime}%
\end{array}
\right)  $, $B=\left(
\begin{array}
[c]{ccc}%
c & c^{\prime} & c^{\prime\prime}\\
d & d^{\prime} & d^{\prime\prime}%
\end{array}
\right)  $, $u=\left(
\begin{array}
[c]{c}%
x\\
y
\end{array}
\right)  $ and $v=\left(
\begin{array}
[c]{c}%
z\\
w
\end{array}
\right)  $, then Proposition \ref{prop.sol.det.rk1upd.claim.AB+uv} says that%
\begin{align*}
&  \left(
\begin{array}
[c]{cccc}%
a & a^{\prime} & a^{\prime\prime} & x\\
b & b^{\prime} & b^{\prime\prime} & y
\end{array}
\right)  \left(
\begin{array}
[c]{cccc}%
c & c^{\prime} & c^{\prime\prime} & z\\
d & d^{\prime} & d^{\prime\prime} & w
\end{array}
\right)  ^{T}\\
&  =\left(
\begin{array}
[c]{ccc}%
a & a^{\prime} & a^{\prime\prime}\\
b & b^{\prime} & b^{\prime\prime}%
\end{array}
\right)  \left(
\begin{array}
[c]{ccc}%
c & c^{\prime} & c^{\prime\prime}\\
d & d^{\prime} & d^{\prime\prime}%
\end{array}
\right)  ^{T}+\left(
\begin{array}
[c]{c}%
x\\
y
\end{array}
\right)  \left(
\begin{array}
[c]{c}%
z\\
w
\end{array}
\right)  ^{T}.
\end{align*}

\end{example}

\begin{proof}
[First proof of Proposition \ref{prop.sol.det.rk1upd.claim.AB+uv}
(sketched).]We can apply Exercise \ref{exe.block2x2.mult} to $n$, $0$, $m$,
$1$, $n$, $0$, $A$, $u$, $0_{0\times m}$, $0_{0\times1}$, $B^{T}$,
$0_{m\times0}$, $v^{T}$ and $0_{1\times0}$ instead of $n$, $n^{\prime}$, $m$,
$m^{\prime}$, $\ell$, $\ell^{\prime}$, $A$, $B$, $C$, $D$, $A^{\prime}$,
$B^{\prime}$, $C^{\prime}$ and $D^{\prime}$. As a result, we obtain%
\[
\left(
\begin{array}
[c]{cc}%
A & u\\
0_{0\times m} & 0_{0\times1}%
\end{array}
\right)  \left(
\begin{array}
[c]{cc}%
B^{T} & 0_{m\times0}\\
v^{T} & 0_{1\times0}%
\end{array}
\right)  =\left(
\begin{array}
[c]{cc}%
AB^{T}+uv^{T} & A0_{m\times0}+u0_{1\times0}\\
0_{0\times m}B^{T}+0_{0\times1}v^{T} & 0_{0\times m}0_{m\times0}+0_{0\times
1}0_{1\times0}%
\end{array}
\right)
\]
(where we are using the notations introduced in Definition \ref{def.block2x2}%
). In view of the equalities%
\[
\left(
\begin{array}
[c]{cc}%
A & u\\
0_{0\times m} & 0_{0\times1}%
\end{array}
\right)  =\left(  A\mid u\right)  \ \ \ \ \ \ \ \ \ \ \left(  \text{this is
rather obvious}\right)  ,
\]%
\[
\left(
\begin{array}
[c]{cc}%
B^{T} & 0_{m\times0}\\
v^{T} & 0_{1\times0}%
\end{array}
\right)  =\left(  B\mid v\right)  ^{T}\ \ \ \ \ \ \ \ \ \ \left(  \text{this
is easy to see}\right)
\]
and%
\begin{align*}
&  \left(
\begin{array}
[c]{cc}%
AB^{T}+uv^{T} & A0_{m\times0}+u0_{1\times0}\\
0_{0\times m}B^{T}+0_{0\times1}v^{T} & 0_{0\times m}0_{m\times0}+0_{0\times
1}0_{1\times0}%
\end{array}
\right) \\
&  =\left(
\begin{array}
[c]{cc}%
AB^{T}+uv^{T} & 0_{n\times0}\\
0_{0\times n} & 0_{0\times0}%
\end{array}
\right)  =AB^{T}+uv^{T},
\end{align*}
this rewrites as $\left(  A\mid u\right)  \left(  B\mid v\right)  ^{T}%
=AB^{T}+uv^{T}$. This proves Proposition \ref{prop.sol.det.rk1upd.claim.AB+uv}.
\end{proof}

We shall now give another, self-contained proof of Proposition
\ref{prop.sol.det.rk1upd.claim.AB+uv}, based upon the following simple fact:

\begin{proposition}
\label{prop.sol.det.rk1upd.claim.AB}Let $n\in\mathbb{N}$ and $m\in\mathbb{N}$.
Let $A\in\mathbb{K}^{n\times m}$ and $B\in\mathbb{K}^{n\times m}$. Then,%
\[
AB^{T}=\sum_{k=1}^{m}A_{\bullet,k}\left(  B_{\bullet,k}\right)  ^{T}.
\]
(Here, we are using the notations introduced in Definition \ref{def.unrows}.)
\end{proposition}

\begin{proof}
[Proof of Proposition \ref{prop.sol.det.rk1upd.claim.AB}.]The matrix $A$ is an
$n\times m$-matrix (since $A\in\mathbb{K}^{n\times m}$). Write this $n\times
m$-matrix $A$ in the form $A=\left(  a_{i,j}\right)  _{1\leq i\leq n,\ 1\leq
j\leq m}$.

The matrix $B$ is an $n\times m$-matrix (since $B\in\mathbb{K}^{n\times m}$).
Write this $n\times m$-matrix $B$ in the form $B=\left(  b_{i,j}\right)
_{1\leq i\leq n,\ 1\leq j\leq m}$. From $B=\left(  b_{i,j}\right)  _{1\leq
i\leq n,\ 1\leq j\leq m}$, we obtain $B^{T}=\left(  b_{j,i}\right)  _{1\leq
i\leq m,\ 1\leq j\leq n}$ (by the definition of $B^{T}$).

The definition of the product of two matrices yields%
\begin{equation}
AB^{T}=\left(  \sum_{k=1}^{m}a_{i,k}b_{j,k}\right)  _{1\leq i\leq n,\ 1\leq
j\leq n} \label{pf.prop.sol.det.rk1upd.claim.AB.0}%
\end{equation}
(since $A=\left(  a_{i,j}\right)  _{1\leq i\leq n,\ 1\leq j\leq m}$ and
$B^{T}=\left(  b_{j,i}\right)  _{1\leq i\leq m,\ 1\leq j\leq n}$).

Let us now fix $p\in\left\{  1,2,\ldots,m\right\}  $. Then, $A_{\bullet,p}$ is
the $p$-th column of the matrix $A$ (by the definition of $A_{\bullet,p}$).
Thus,%
\begin{align*}
A_{\bullet,p}  &  =\left(  \text{the }p\text{-th column of the matrix
}A\right) \\
&  =\left(
\begin{array}
[c]{c}%
a_{1,p}\\
a_{2,p}\\
\vdots\\
a_{n,p}%
\end{array}
\right)  \ \ \ \ \ \ \ \ \ \ \left(  \text{since }A=\left(  a_{i,j}\right)
_{1\leq i\leq n,\ 1\leq j\leq m}\right) \\
&  =\left(  a_{i,p}\right)  _{1\leq i\leq n,\ 1\leq j\leq1}.
\end{align*}
The same argument (applied to $B$ and $b_{i,j}$ instead of $A$ and $a_{i,j}$)
yields $B_{\bullet,p}=\left(  b_{i,p}\right)  _{1\leq i\leq n,\ 1\leq j\leq1}%
$. Hence, $\left(  B_{\bullet,p}\right)  ^{T}=\left(  b_{j,p}\right)  _{1\leq
i\leq1,\ 1\leq j\leq n}$ (by the definition of $\left(  B_{\bullet,p}\right)
^{T}$).

Now, the definition of the product of two matrices yields%
\begin{align}
A_{\bullet,p}\left(  B_{\bullet,p}\right)  ^{T}  &  =\left(  \underbrace{\sum
_{k=1}^{1}a_{i,p}b_{j,p}}_{=a_{i,p}b_{j,p}}\right)  _{1\leq i\leq n,\ 1\leq
j\leq n}\ \ \ \ \ \ \ \ \ \ \left(
\begin{array}
[c]{c}%
\text{since }A_{\bullet,p}=\left(  a_{i,p}\right)  _{1\leq i\leq n,\ 1\leq
j\leq1}\\
\text{and }\left(  B_{\bullet,p}\right)  ^{T}=\left(  b_{j,p}\right)  _{1\leq
i\leq1,\ 1\leq j\leq n}%
\end{array}
\right) \nonumber\\
&  =\left(  a_{i,p}b_{j,p}\right)  _{1\leq i\leq n,\ 1\leq j\leq n}.
\label{pf.prop.sol.det.rk1upd.claim.AB.1}%
\end{align}

Now, let us forget that we fixed $p$. We thus have proven
(\ref{pf.prop.sol.det.rk1upd.claim.AB.1}) for each $p\in\left\{
1,2,\ldots,m\right\}  $. Now,%
\[
\sum_{k=1}^{m}\underbrace{A_{\bullet,k}\left(  B_{\bullet,k}\right)  ^{T}%
}_{\substack{=\left(  a_{i,k}b_{j,k}\right)  _{1\leq i\leq n,\ 1\leq j\leq
n}\\\text{(by (\ref{pf.prop.sol.det.rk1upd.claim.AB.1}), applied to
}p=k\text{)}}}=\sum_{k=1}^{m}\left(  a_{i,k}b_{j,k}\right)  _{1\leq i\leq
n,\ 1\leq j\leq n}=\left(  \sum_{k=1}^{m}a_{i,k}b_{j,k}\right)  _{1\leq i\leq
n,\ 1\leq j\leq n}.
\]
Comparing this with (\ref{pf.prop.sol.det.rk1upd.claim.AB.0}), we obtain
$AB^{T}=\sum_{k=1}^{m}A_{\bullet,k}\left(  B_{\bullet,k}\right)  ^{T}$. This
proves Proposition \ref{prop.sol.det.rk1upd.claim.AB}.
\end{proof}

\begin{proof}
[Second proof of Proposition \ref{prop.sol.det.rk1upd.claim.AB+uv}.]We have
$\left(  A\mid u\right)  \in\mathbb{K}^{n\times\left(  m+1\right)  }$ (since
$A\in\mathbb{K}^{n\times m}$ and $u\in\mathbb{K}^{n\times1}$) and $\left(
B\mid v\right)  \in\mathbb{K}^{n\times\left(  m+1\right)  }$ (since
$B\in\mathbb{K}^{n\times m}$ and $v\in\mathbb{K}^{n\times1}$). Thus,
Proposition \ref{prop.sol.det.rk1upd.claim.AB} (applied to $m+1$, $\left(
A\mid u\right)  $ and $\left(  B\mid v\right)  $ instead of $m$, $A$ and $B$)
yields%
\begin{align*}
&  \left(  A\mid u\right)  \left(  B\mid v\right)  ^{T}\\
&  =\sum_{k=1}^{m+1}\left(  A\mid u\right)  _{\bullet,k}\left(  \left(  B\mid
v\right)  _{\bullet,k}\right)  ^{T}\\
&  =\sum_{k=1}^{m}\underbrace{\left(  A\mid u\right)  _{\bullet,k}%
}_{\substack{=A_{\bullet,k}\\\text{(by Proposition \ref{prop.addcol.props1}
\textbf{(a)},}\\\text{applied to }u\text{ and }k\text{ instead of }v\text{ and
}q\text{)}}}\left(  \underbrace{\left(  B\mid v\right)  _{\bullet,k}%
}_{\substack{=B_{\bullet,k}\\\text{(by Proposition \ref{prop.addcol.props1}
\textbf{(a)},}\\\text{applied to }B\text{ and }k\text{ instead of }A\text{ and
}q\text{)}}}\right)  ^{T}\\
&  \ \ \ \ \ \ \ \ \ \ +\underbrace{\left(  A\mid u\right)  _{\bullet,m+1}%
}_{\substack{=u\\\text{(by Proposition \ref{prop.addcol.props1} \textbf{(b)}%
,}\\\text{applied to }u\text{ instead of }v\text{)}}}\left(
\underbrace{\left(  B\mid v\right)  _{\bullet,m+1}}_{\substack{=v\\\text{(by
Proposition \ref{prop.addcol.props1} \textbf{(b)},}\\\text{applied to }B\text{
instead of }A\text{)}}}\right)  ^{T}\\
&  \ \ \ \ \ \ \ \ \ \ \left(  \text{here, we have split off the addend for
}k=m+1\text{ from the sum}\right) \\
&  =\sum_{k=1}^{m}A_{\bullet,k}\left(  B_{\bullet,k}\right)  ^{T}+uv^{T}.
\end{align*}
Comparing this with%
\[
\underbrace{AB^{T}}_{\substack{=\sum_{k=1}^{m}A_{\bullet,k}\left(
B_{\bullet,k}\right)  ^{T}\\\text{(by Proposition
\ref{prop.sol.det.rk1upd.claim.AB})}}}+uv^{T}=\sum_{k=1}^{m}A_{\bullet
,k}\left(  B_{\bullet,k}\right)  ^{T}+uv^{T},
\]
we obtain $\left(  A\mid u\right)  \left(  B\mid v\right)  ^{T}=AB^{T}+uv^{T}%
$. Thus, Proposition \ref{prop.sol.det.rk1upd.claim.AB+uv} is proven again.
\end{proof}

\begin{corollary}
\label{cor.sol.det.rk1upd.claim.A+uv}Let $n\in\mathbb{N}$. Let $u\in
\mathbb{K}^{1\times n}$ and $v\in\mathbb{K}^{n\times1}$. Let $A$ be an
$n\times n$-matrix. Then,%
\[
A+uv=\left(  A\mid u\right)  \left(  I_{n}\mid v^{T}\right)  ^{T}.
\]
(Here, we are using the notations introduced in Definition \ref{def.addcol}.)
\end{corollary}

\begin{proof}
[Proof of Corollary \ref{cor.sol.det.rk1upd.claim.A+uv}.]We have
$v\in\mathbb{K}^{n\times1}$ and thus $v^{T}\in\mathbb{K}^{1\times n}$. Also,
$I_{n}\in\mathbb{K}^{n\times n}$. Thus, Proposition
\ref{prop.sol.det.rk1upd.claim.AB+uv} (applied to $n$, $I_{n}$ and $v^{T}$
instead of $m$, $B$ and $v$) yields%
\[
\left(  A\mid u\right)  \left(  I_{n}\mid v^{T}\right)  ^{T}%
=A\underbrace{\left(  I_{n}\right)  ^{T}}_{=I_{n}}+u\underbrace{\left(
v^{T}\right)  ^{T}}_{\substack{=v\\\text{(since }\left(  C^{T}\right)
^{T}=C\text{ for any}\\\text{matrix }C\text{)}}}=\underbrace{AI_{n}}%
_{=A}+uv=A+uv.
\]
This proves Corollary \ref{cor.sol.det.rk1upd.claim.A+uv}.
\end{proof}

Next, we state some simple facts about determinants:

\begin{lemma}
\label{lem.sol.det.rk1upd.adjI}Let $n\in\mathbb{N}$. Then:

\textbf{(a)} Every $u\in\left\{  1,2,\ldots,n\right\}  $ satisfies
$\det\left(  \left(  I_{n}\right)  _{\sim u,\sim u}\right)  =1$.

\textbf{(b)} If $u$ and $v$ are two elements of $\left\{  1,2,\ldots
,n\right\}  $ such that $u\neq v$, then $\det\left(  \left(  I_{n}\right)
_{\sim u,\sim v}\right)  =0$.
\end{lemma}

\begin{proof}
[Proof of Lemma \ref{lem.sol.det.rk1upd.adjI}.]There are myriad ways to prove
Lemma \ref{lem.sol.det.rk1upd.adjI} (for example, one can prove part
\textbf{(a)} by showing that $\left(  I_{n}\right)  _{\sim u,\sim u}=I_{n-1}$,
and prove part \textbf{(b)} by arguing that the matrix $\left(  I_{n}\right)
_{\sim u,\sim v}$ has a row consisting of zeroes\footnote{or, what also
suffices, a column consisting of zeroes}). The proof we will now give is not
the simplest one, but the shortest one (using what we have proven so far):

For any two objects $i$ and $j$, define an element $\delta_{i,j}\in\mathbb{K}$
by $\delta_{i,j}=%
\begin{cases}
1, & \text{if }i=j;\\
0, & \text{if }i\neq j
\end{cases}
$. Then, $I_{n}=\left(  \delta_{i,j}\right)  _{1\leq i\leq n,\ 1\leq j\leq n}$
(by the definition of $I_{n}$).

Now, Theorem \ref{thm.adj.inverse} (applied to $A=I_{n}$) yields $I_{n}%
\cdot\operatorname*{adj}\left(  I_{n}\right)  =\operatorname*{adj}\left(
I_{n}\right)  \cdot I_{n}=\det\left(  I_{n}\right)  \cdot I_{n}$. Thus,
$\operatorname*{adj}\left(  I_{n}\right)  \cdot I_{n}=\underbrace{\det\left(
I_{n}\right)  }_{=1}\cdot I_{n}=I_{n}$, so that $I_{n}=\operatorname*{adj}%
\left(  I_{n}\right)  \cdot I_{n}=\operatorname*{adj}\left(  I_{n}\right)  $
(since $BI_{n}=I_{n}$ for any $m\in\mathbb{N}$ and any $m\times n$-matrix
$B$). Thus,%
\[
I_{n}=\operatorname*{adj}\left(  I_{n}\right)  =\left(  \left(  -1\right)
^{i+j}\det\left(  \left(  I_{n}\right)  _{\sim j,\sim i}\right)  \right)
_{1\leq i\leq n,\ 1\leq j\leq n}%
\]
(by the definition of $\operatorname*{adj}\left(  I_{n}\right)  $). Hence,%
\[
\left(  \left(  -1\right)  ^{i+j}\det\left(  \left(  I_{n}\right)  _{\sim
j,\sim i}\right)  \right)  _{1\leq i\leq n,\ 1\leq j\leq n}=I_{n}=\left(
\delta_{i,j}\right)  _{1\leq i\leq n,\ 1\leq j\leq n}.
\]
In other words,%
\begin{equation}
\left(  -1\right)  ^{i+j}\det\left(  \left(  I_{n}\right)  _{\sim j,\sim
i}\right)  =\delta_{i,j} \label{pf.lem.sol.det.rk1upd.adjI.1}%
\end{equation}
for every $i\in\left\{  1,2,\ldots,n\right\}  $ and $j\in\left\{
1,2,\ldots,n\right\}  $.

\textbf{(a)} Let $u\in\left\{  1,2,\ldots,n\right\}  $. Then,
(\ref{pf.lem.sol.det.rk1upd.adjI.1}) (applied to $i=u$ and $j=u$) yields%
\[
\left(  -1\right)  ^{u+u}\det\left(  \left(  I_{n}\right)  _{\sim u,\sim
u}\right)  =\delta_{u,u}=1\ \ \ \ \ \ \ \ \ \ \left(  \text{since }u=u\right)
.
\]
Comparing this with $\underbrace{\left(  -1\right)  ^{u+u}}%
_{\substack{=1\\\text{(since }u+u=2u\text{ is even)}}}\det\left(  \left(
I_{n}\right)  _{\sim u,\sim u}\right)  =\det\left(  \left(  I_{n}\right)
_{\sim u,\sim u}\right)  $, we obtain $\det\left(  \left(  I_{n}\right)
_{\sim u,\sim u}\right)  =1$. This proves Lemma \ref{lem.sol.det.rk1upd.adjI}
\textbf{(a)}.

\textbf{(b)} Let $u$ and $v$ be two elements of $\left\{  1,2,\ldots
,n\right\}  $ such that $u\neq v$. Then, (\ref{pf.lem.sol.det.rk1upd.adjI.1})
(applied to $i=v$ and $j=u$) yields%
\[
\left(  -1\right)  ^{v+u}\det\left(  \left(  I_{n}\right)  _{\sim u,\sim
v}\right)  =\delta_{v,u}=0\ \ \ \ \ \ \ \ \ \ \left(  \text{since }v\neq
u\text{ (since }u\neq v\text{)}\right)  .
\]
Multiplying both sides of this equality by $\left(  -1\right)  ^{v+u}$, we
obtain%
\[
\left(  -1\right)  ^{v+u}\left(  -1\right)  ^{v+u}\det\left(  \left(
I_{n}\right)  _{\sim u,\sim v}\right)  =0.
\]
Comparing this with
\[
\underbrace{\left(  -1\right)  ^{v+u}\left(  -1\right)  ^{v+u}}%
_{\substack{=\left(  -1\right)  ^{\left(  v+u\right)  +\left(  v+u\right)
}=1\\\text{(since }\left(  v+u\right)  +\left(  v+u\right)  =2\left(
v+u\right)  \text{ is even)}}}\det\left(  \left(  I_{n}\right)  _{\sim u,\sim
v}\right)  =\det\left(  \left(  I_{n}\right)  _{\sim u,\sim v}\right)  ,
\]
we obtain $\det\left(  \left(  I_{n}\right)  _{\sim u,\sim v}\right)  =0$.
Lemma \ref{lem.sol.det.rk1upd.adjI} \textbf{(b)} is thus proven.
\end{proof}

\begin{proposition}
\label{prop.sol.det.rk1upd.Iv}Let $n\in\mathbb{N}$. Let $v=\left(  v_{1}%
,v_{2},\ldots,v_{n}\right)  ^{T}\in\mathbb{K}^{n\times1}$ and $k\in\left\{
1,2,\ldots,n\right\}  $. Then,%
\[
\det\left(  \left(  I_{n}\mid v\right)  _{\bullet,\sim k}\right)  =\left(
-1\right)  ^{n+k}v_{k}.
\]
(Here, we are using the notations introduced in Definition \ref{def.unrows}.)
\end{proposition}

\begin{proof}
[Proof of Proposition \ref{prop.sol.det.rk1upd.Iv}.]We have $I_{n}%
\in\mathbb{K}^{n\times n}$ and thus $\left(  I_{n}\right)  _{\bullet,\sim
k}\in\mathbb{K}^{n\times\left(  n-1\right)  }$. Every $i\in\left\{
1,2,\ldots,n\right\}  $ satisfies%
\begin{equation}
\left(  \left(  I_{n}\right)  _{\bullet,\sim k}\right)  _{\sim i,\bullet
}=\left(  I_{n}\right)  _{\sim i,\sim k} \label{pf.prop.sol.det.rk1upd.Iv.1}%
\end{equation}
\footnote{\textit{Proof of (\ref{pf.prop.sol.det.rk1upd.Iv.1}):} Let
$i\in\left\{  1,2,\ldots,n\right\}  $. Proposition \ref{prop.unrows.basics}
\textbf{(c)} (applied to $n$, $I_{n}$, $i$ and $k$ instead of $m$, $A$, $u$
and $v$) shows that $\left(  \left(  I_{n}\right)  _{\bullet,\sim k}\right)
_{\sim i,\bullet}=\left(  \left(  I_{n}\right)  _{\sim i,\bullet}\right)
_{\bullet,\sim k}=\left(  I_{n}\right)  _{\sim i,\sim k}$. This proves
(\ref{pf.prop.sol.det.rk1upd.Iv.1}).}.

Proposition \ref{prop.addcol.props1} \textbf{(d)} (applied to $m=n$, $A=I_{n}$
and $q=k$) yields $\left(  I_{n}\mid v\right)  _{\bullet,\sim k}=\left(
\left(  I_{n}\right)  _{\bullet,\sim k}\mid v\right)  $. Hence,%
\begin{align*}
&  \det\underbrace{\left(  \left(  I_{n}\mid v\right)  _{\bullet,\sim
k}\right)  }_{=\left(  \left(  I_{n}\right)  _{\bullet,\sim k}\mid v\right)
}\\
&  =\det\left(  \left(  I_{n}\right)  _{\bullet,\sim k}\mid v\right) \\
&  =\underbrace{\sum_{i=1}^{n}}_{=\sum_{i\in\left\{  1,2,\ldots,n\right\}  }%
}\left(  -1\right)  ^{n+i}v_{i}\det\left(  \underbrace{\left(  \left(
I_{n}\right)  _{\bullet,\sim k}\right)  _{\sim i,\bullet}}_{\substack{=\left(
I_{n}\right)  _{\sim i,\sim k}\\\text{(by (\ref{pf.prop.sol.det.rk1upd.Iv.1}%
))}}}\right) \\
&  \ \ \ \ \ \ \ \ \ \ \left(  \text{by Proposition \ref{prop.addcol.props2}
\textbf{(a)}, applied to }A=\left(  I_{n}\right)  _{\bullet,\sim k}\right) \\
&  =\sum_{i\in\left\{  1,2,\ldots,n\right\}  }\left(  -1\right)  ^{n+i}%
v_{i}\det\left(  \left(  I_{n}\right)  _{\sim i,\sim k}\right) \\
&  =\left(  -1\right)  ^{n+k}v_{k}\underbrace{\det\left(  \left(
I_{n}\right)  _{\sim k,\sim k}\right)  }_{\substack{=1\\\text{(by Lemma
\ref{lem.sol.det.rk1upd.adjI} \textbf{(a)},}\\\text{applied to }k\text{
instead of }u\text{)}}}+\sum_{\substack{i\in\left\{  1,2,\ldots,n\right\}
;\\i\neq k}}\left(  -1\right)  ^{n+i}v_{i}\underbrace{\det\left(  \left(
I_{n}\right)  _{\sim i,\sim k}\right)  }_{\substack{=0\\\text{(by Lemma
\ref{lem.sol.det.rk1upd.adjI} \textbf{(b)},}\\\text{applied to }i\text{ and
}k\\\text{instead of }u\text{ and }v\text{)}}}\\
&  \ \ \ \ \ \ \ \ \ \ \left(  \text{here, we have split off the addend for
}i=k\text{ from the sum}\right) \\
&  =\left(  -1\right)  ^{n+k}v_{k}+\underbrace{\sum_{\substack{i\in\left\{
1,2,\ldots,n\right\}  ;\\i\neq k}}\left(  -1\right)  ^{n+i}v_{i}0}%
_{=0}=\left(  -1\right)  ^{n+k}v_{k}.
\end{align*}
This proves Proposition \ref{prop.sol.det.rk1upd.Iv}.
\end{proof}

\begin{proposition}
\label{prop.sol.det.rk1upd.det(ABT)}Let $n\in\mathbb{N}$. Let $A\in
\mathbb{K}^{n\times\left(  n+1\right)  }$ and $B\in\mathbb{K}^{n\times\left(
n+1\right)  }$. Then,%
\[
\det\left(  AB^{T}\right)  =\sum_{k=1}^{n+1}\det\left(  A_{\bullet,\sim
k}\right)  \det\left(  B_{\bullet,\sim k}\right)  .
\]
(Here, we are using the notations introduced in Definition \ref{def.unrows}.)
\end{proposition}

\begin{proof}
[Proof of Proposition \ref{prop.sol.det.rk1upd.det(ABT)}.]Let us use the
notations introduced in Definition \ref{def.rowscols}. Clearly, $n+1$ is a
positive integer (since $n\in\mathbb{N}$).

Also, $A$ is an $n\times\left(  n+1\right)  $-matrix (since $A\in
\mathbb{K}^{n\times\left(  n+1\right)  }$). In other words, $A$ is an $\left(
\left(  n+1\right)  -1\right)  \times\left(  n+1\right)  $-matrix (since
$n=\left(  n+1\right)  -1$).

Also, $B$ is an $n\times\left(  n+1\right)  $-matrix (since $B\in
\mathbb{K}^{n\times\left(  n+1\right)  }$). Hence, $B^{T}$ is an $\left(
n+1\right)  \times n$-matrix. In other words, $B^{T}$ is an $\left(
n+1\right)  \times\left(  \left(  n+1\right)  -1\right)  $-matrix (since
$n=\left(  n+1\right)  -1$).

Now, every $k\in\left\{  1,2,\ldots,n+1\right\}  $ satisfies%
\begin{equation}
\operatorname*{cols}\nolimits_{1,2,\ldots,\widehat{k},\ldots,n+1}%
A=A_{\bullet,\sim k} \label{pf.prop.sol.det.rk1upd.det(ABT).1}%
\end{equation}
\footnote{\textit{Proof of (\ref{pf.prop.sol.det.rk1upd.det(ABT).1}):} Let
$k\in\left\{  1,2,\ldots,n+1\right\}  $. Then, the definition of
$A_{\bullet,\sim k}$ yields $A_{\bullet,\sim k}=\operatorname*{cols}%
\nolimits_{1,2,\ldots,\widehat{k},\ldots,n+1}A$ (since $A\in\mathbb{K}%
^{n\times\left(  n+1\right)  }$). This proves
(\ref{pf.prop.sol.det.rk1upd.det(ABT).1}).} and%
\begin{equation}
\operatorname*{rows}\nolimits_{1,2,\ldots,\widehat{k},\ldots,n+1}\left(
B^{T}\right)  =\left(  B_{\bullet,\sim k}\right)  ^{T}
\label{pf.prop.sol.det.rk1upd.det(ABT).2}%
\end{equation}
\footnote{\textit{Proof of (\ref{pf.prop.sol.det.rk1upd.det(ABT).2}):} Let
$k\in\left\{  1,2,\ldots,n+1\right\}  $. Then, the definition of $\left(
B^{T}\right)  _{\sim k,\bullet}$ yields $\left(  B^{T}\right)  _{\sim
k,\bullet}=\operatorname*{rows}\nolimits_{1,2,\ldots,\widehat{k},\ldots
,n+1}\left(  B^{T}\right)  $ (since $B^{T}\in\mathbb{K}^{\left(  n+1\right)
\times n}$ (since $B^{T}$ is an $\left(  n+1\right)  \times n$-matrix)). But
Lemma \ref{lem.unrows.transpose.1} (applied to $m=n+1$ and $r=k$) yields
$\left(  B^{T}\right)  _{\sim k,\bullet}=\left(  B_{\bullet,\sim k}\right)
^{T}$. Comparing this with $\left(  B^{T}\right)  _{\sim k,\bullet
}=\operatorname*{rows}\nolimits_{1,2,\ldots,\widehat{k},\ldots,n+1}\left(
B^{T}\right)  $, we obtain $\operatorname*{rows}\nolimits_{1,2,\ldots
,\widehat{k},\ldots,n+1}\left(  B^{T}\right)  =\left(  B_{\bullet,\sim
k}\right)  ^{T}$. This proves (\ref{pf.prop.sol.det.rk1upd.det(ABT).2}).} and%
\begin{equation}
\det\left(  \operatorname*{rows}\nolimits_{1,2,\ldots,\widehat{k},\ldots
,n+1}\left(  B^{T}\right)  \right)  =\det\left(  B_{\bullet,\sim k}\right)
\label{pf.prop.sol.det.rk1upd.det(ABT).3}%
\end{equation}
\footnote{\textit{Proof of (\ref{pf.prop.sol.det.rk1upd.det(ABT).3}):} Let
$k\in\left\{  1,2,\ldots,n+1\right\}  $. From $B\in\mathbb{K}^{n\times\left(
n+1\right)  }$, we obtain $B_{\bullet,\sim k}\in\mathbb{K}^{n\times n}$. In
other words, $B_{\bullet,\sim k}$ is an $n\times n$-matrix. Hence, Exercise
\ref{exe.ps4.4} (applied to $B_{\bullet,\sim k}$ instead of $A$) yields
$\det\left(  \left(  B_{\bullet,\sim k}\right)  ^{T}\right)  =\det\left(
B_{\bullet,\sim k}\right)  $. Now,%
\[
\det\left(  \underbrace{\operatorname*{rows}\nolimits_{1,2,\ldots
,\widehat{k},\ldots,n+1}\left(  B^{T}\right)  }_{\substack{=\left(
B_{\bullet,\sim k}\right)  ^{T}\\\text{(by
(\ref{pf.prop.sol.det.rk1upd.det(ABT).2}))}}}\right)  =\det\left(  \left(
B_{\bullet,\sim k}\right)  ^{T}\right)  =\det\left(  B_{\bullet,\sim
k}\right)  .
\]
This proves (\ref{pf.prop.sol.det.rk1upd.det(ABT).3}).}.

Now, Lemma \ref{lem.adj(AB).cauchy-binet} (applied to $n+1$ and $B^{T}$
instead of $n$ and $B$) yields%
\begin{align*}
\det\left(  AB^{T}\right)   &  =\sum_{k=1}^{n+1}\det\left(
\underbrace{\operatorname*{cols}\nolimits_{1,2,\ldots,\widehat{k},\ldots
,n+1}A}_{\substack{=A_{\bullet,\sim k}\\\text{(by
(\ref{pf.prop.sol.det.rk1upd.det(ABT).1}))}}}\right)  \cdot\underbrace{\det
\left(  \operatorname*{rows}\nolimits_{1,2,\ldots,\widehat{k},\ldots
,n+1}\left(  B^{T}\right)  \right)  }_{\substack{=\det\left(  B_{\bullet,\sim
k}\right)  \\\text{(by (\ref{pf.prop.sol.det.rk1upd.det(ABT).3}))}}}\\
&  =\sum_{k=1}^{n+1}\det\left(  A_{\bullet,\sim k}\right)  \det\left(
B_{\bullet,\sim k}\right)  .
\end{align*}
This proves Proposition \ref{prop.sol.det.rk1upd.det(ABT)}.
\end{proof}

The next lemma is an easy consequence of results proven before:

\begin{lemma}
\label{lem.sol.det.rk1upd.detAuk}Let $n\in\mathbb{N}$. Let $u=\left(
u_{1},u_{2},\ldots,u_{n}\right)  ^{T}\in\mathbb{K}^{n\times1}$ be a column
vector with $n$ entries. Let $A$ be an $n\times n$-matrix. Let $k\in\left\{
1,2,\ldots,n\right\}  $. Then,%
\[
\det\left(  \left(  A\mid u\right)  _{\bullet,\sim k}\right)  =\sum_{j=1}%
^{n}\left(  -1\right)  ^{n+j}u_{j}\det\left(  A_{\sim j,\sim k}\right)  .
\]

\end{lemma}

\begin{proof}
[Proof of Lemma \ref{lem.sol.det.rk1upd.detAuk}.]From $A\in\mathbb{K}^{n\times
n}$ and $u\in\mathbb{K}^{n\times1}$, we obtain $\left(  A\mid u\right)
\in\mathbb{K}^{n\times\left(  n+1\right)  }$.

From $k\in\left\{  1,2,\ldots,n\right\}  $, we obtain $A_{\bullet,\sim k}%
\in\mathbb{K}^{n\times\left(  n-1\right)  }$ (since $A\in\mathbb{K}^{n\times
n}$). Proposition \ref{prop.addcol.props1} \textbf{(c)} (applied to $n$, $u$
and $k$ instead of $m$, $v$ and $q$) shows that $\left(  A\mid u\right)
_{\bullet,\sim k}=\left(  A_{\bullet,\sim k}\mid u\right)  $. Hence,%
\begin{align*}
\det\left(  \underbrace{\left(  A\mid u\right)  _{\bullet,\sim k}}_{=\left(
A_{\bullet,\sim k}\mid u\right)  }\right)   &  =\det\left(  A_{\bullet,\sim
k}\mid u\right) \\
&  =\sum_{i=1}^{n}\left(  -1\right)  ^{n+i}u_{i}\det\left(
\underbrace{\left(  A_{\bullet,\sim k}\right)  _{\sim i,\bullet}%
}_{\substack{=A_{\sim i,\sim k}\\\text{(since Proposition
\ref{prop.unrows.basics} \textbf{(c)} (applied to}\\n\text{, }i\text{ and
}k\text{ instead of }m\text{, }u\text{ and }v\text{) shows that}\\\left(
A_{\bullet,\sim k}\right)  _{\sim i,\bullet}=\left(  A_{\sim i,\bullet
}\right)  _{\bullet,\sim k}=A_{\sim i,\sim k}\text{)}}}\right) \\
&  \ \ \ \ \ \ \ \ \ \ \left(
\begin{array}
[c]{c}%
\text{by Proposition \ref{prop.addcol.props2} \textbf{(a)}, applied to}\\
A_{\bullet,\sim k}\text{, }u\text{ and }u_{i}\text{ instead of }A\text{,
}v\text{ and }v_{i}%
\end{array}
\right) \\
&  =\sum_{i=1}^{n}\left(  -1\right)  ^{n+i}u_{i}\det\left(  A_{\sim i,\sim
k}\right)  =\sum_{j=1}^{n}\left(  -1\right)  ^{n+j}u_{j}\det\left(  A_{\sim
j,\sim k}\right)
\end{align*}
(here, we have renamed the summation index $i$ as $j$). This proves Lemma
\ref{lem.sol.det.rk1upd.detAuk}.
\end{proof}

Now we are more than ready to easily prove Theorem
\ref{thm.sol.det.rk1upd.claim}:

\begin{proof}
[Proof of Theorem \ref{thm.sol.det.rk1upd.claim}.]We shall use the notations
introduced in Definition \ref{def.submatrix.minor}, in Definition
\ref{def.addcol} and in Definition \ref{def.unrows}.

We know that $u$ is a column vector with $n$ entries; in other words,
$u\in\mathbb{K}^{n\times1}$. Also, $v$ is a row vector with $n$ entries; in
other words, $v\in\mathbb{K}^{1\times n}$. Hence, $v^{T}\in\mathbb{K}%
^{n\times1}$.

From $A\in\mathbb{K}^{n\times n}$ and $u\in\mathbb{K}^{n\times1}$, we obtain
$\left(  A\mid u\right)  \in\mathbb{K}^{n\times\left(  n+1\right)  }$. From
$I_{n}\in\mathbb{K}^{n\times n}$ and $v^{T}\in\mathbb{K}^{n\times1}$, we
obtain $\left(  I_{n}\mid v^{T}\right)  \in\mathbb{K}^{n\times\left(
n+1\right)  }$.

Write the vector $v$ in the form $v=\left(  v_{1},v_{2},\ldots,v_{n}\right)
$. (This is possible, since $v$ is a row vector with $n$ entries.) From
$v=\left(  v_{1},v_{2},\ldots,v_{n}\right)  $, we obtain $v^{T}=\left(
v_{1},v_{2},\ldots,v_{n}\right)  ^{T}$.

Write the vector $u$ in the form $u=\left(  u_{1},u_{2},\ldots,u_{n}\right)
^{T}$. (This is possible, since $u$ is a column vector with $n$ entries.)

Now,%
\begin{align*}
&  \det\left(  \underbrace{A+uv}_{\substack{=\left(  A\mid u\right)  \left(
I_{n}\mid v^{T}\right)  ^{T}\\\text{(by Corollary
\ref{cor.sol.det.rk1upd.claim.A+uv})}}}\right) \\
&  =\det\left(  \left(  A\mid u\right)  \left(  I_{n}\mid v^{T}\right)
^{T}\right)  =\sum_{k=1}^{n+1}\det\left(  \left(  A\mid u\right)
_{\bullet,\sim k}\right)  \det\left(  \left(  I_{n}\mid v^{T}\right)
_{\bullet,\sim k}\right) \\
&  \ \ \ \ \ \ \ \ \ \ \left(
\begin{array}
[c]{c}%
\text{by Proposition \ref{prop.sol.det.rk1upd.det(ABT)}, applied to }\left(
A\mid u\right)  \text{ and }\left(  I_{n}\mid v^{T}\right) \\
\text{instead of }A\text{ and }B
\end{array}
\right) \\
&  =\sum_{k=1}^{n}\underbrace{\det\left(  \left(  A\mid u\right)
_{\bullet,\sim k}\right)  }_{\substack{=\sum_{j=1}^{n}\left(  -1\right)
^{n+j}u_{j}\det\left(  A_{\sim j,\sim k}\right)  \\\text{(by Lemma
\ref{lem.sol.det.rk1upd.detAuk})}}}\underbrace{\det\left(  \left(  I_{n}\mid
v^{T}\right)  _{\bullet,\sim k}\right)  }_{\substack{=\left(  -1\right)
^{n+k}v_{k}\\\text{(by Proposition \ref{prop.sol.det.rk1upd.Iv},}%
\\\text{applied to }v^{T}\text{ instead of }v\text{)}}}\\
&  \ \ \ \ \ \ \ \ \ \ +\det\left(  \underbrace{\left(  A\mid u\right)
_{\bullet,\sim\left(  n+1\right)  }}_{\substack{=A\\\text{(by Proposition
\ref{prop.addcol.props1} \textbf{(d)}, applied}\\\text{to }n\text{ and
}u\text{ instead of }m\text{ and }v\text{)}}}\right)  \det\left(
\underbrace{\left(  I_{n}\mid v^{T}\right)  _{\bullet,\sim\left(  n+1\right)
}}_{\substack{=I_{n}\\\text{(by Proposition \ref{prop.addcol.props1}
\textbf{(d)}, applied}\\\text{to }n\text{, }I_{n}\text{ and }v^{T}\text{
instead of }m\text{, }A\text{ and }v\text{)}}}\right) \\
&  \ \ \ \ \ \ \ \ \ \ \left(  \text{here, we have split off the addend for
}k=n+1\text{ from the sum}\right) \\
&  =\underbrace{\sum_{k=1}^{n}\left(  \sum_{j=1}^{n}\left(  -1\right)
^{n+j}u_{j}\det\left(  A_{\sim j,\sim k}\right)  \right)  \cdot\left(
-1\right)  ^{n+k}v_{k}}_{=\sum_{k=1}^{n}\sum_{j=1}^{n}\left(  -1\right)
^{n+j}u_{j}\det\left(  A_{\sim j,\sim k}\right)  \cdot\left(  -1\right)
^{n+k}v_{k}}+\det A\cdot\underbrace{\det\left(  I_{n}\right)  }_{=1}\\
&  =\sum_{k=1}^{n}\sum_{j=1}^{n}\left(  -1\right)  ^{n+j}u_{j}\det\left(
A_{\sim j,\sim k}\right)  \cdot\left(  -1\right)  ^{n+k}v_{k}+\det A.
\end{align*}
Subtracting $\det A$ from both sides of this equality, we obtain%
\begin{align*}
&  \det\left(  A+uv\right)  -\det A\\
&  =\sum_{k=1}^{n}\sum_{j=1}^{n}\left(  -1\right)  ^{n+j}u_{j}\det\left(
A_{\sim j,\sim k}\right)  \cdot\left(  -1\right)  ^{n+k}v_{k}\\
&  =\sum_{i=1}^{n}\sum_{j=1}^{n}\underbrace{\left(  -1\right)  ^{n+j}u_{j}%
\det\left(  A_{\sim j,\sim i}\right)  \cdot\left(  -1\right)  ^{n+i}v_{i}%
}_{=\left(  -1\right)  ^{n+i}\left(  -1\right)  ^{n+j}u_{j}v_{i}\det\left(
A_{\sim j,\sim i}\right)  }\\
&  \ \ \ \ \ \ \ \ \ \ \left(  \text{here, we have renamed the summation index
}k\text{ as }i\text{ in the first sum}\right) \\
&  =\sum_{i=1}^{n}\sum_{j=1}^{n}\underbrace{\left(  -1\right)  ^{n+i}\left(
-1\right)  ^{n+j}}_{\substack{=\left(  -1\right)  ^{\left(  n+i\right)
+\left(  n+j\right)  }=\left(  -1\right)  ^{i+j}\\\text{(since }\left(
n+i\right)  +\left(  n+j\right)  =2n+i+j\equiv i+j\operatorname{mod}2\text{)}%
}}u_{j}v_{i}\det\left(  A_{\sim j,\sim i}\right) \\
&  =\sum_{i=1}^{n}\sum_{j=1}^{n}\left(  -1\right)  ^{i+j}u_{j}v_{i}\det\left(
A_{\sim j,\sim i}\right)  .
\end{align*}
Comparing this with%
\[
\operatorname*{ent}\left(  v\left(  \operatorname*{adj}A\right)  u\right)
=\sum_{i=1}^{n}\sum_{j=1}^{n}\left(  -1\right)  ^{i+j}u_{j}v_{i}\det\left(
A_{\sim j,\sim i}\right)  \ \ \ \ \ \ \ \ \ \ \left(  \text{by Corollary
\ref{cor.sol.det.rk1upd.claim.ent(v(adjA)u)}}\right)  ,
\]
we obtain $\det\left(  A+uv\right)  -\det A=\operatorname*{ent}\left(
v\left(  \operatorname*{adj}A\right)  u\right)  $. In other words,
$\det\left(  A+uv\right)  =\det A+\operatorname*{ent}\left(  v\left(
\operatorname*{adj}A\right)  u\right)  $. This proves Theorem
\ref{thm.sol.det.rk1upd.claim}.
\end{proof}

Now, Theorem \ref{thm.sol.det.rk1upd.claim} is proven; in other words,
Exercise \ref{exe.det.rk1upd} is solved (since Theorem
\ref{thm.sol.det.rk1upd.claim} is just a restatement of the claim of this exercise).

\subsection{Solution to Exercise \ref{exe.det.bordered}}

In this section, we shall use the notations introduced in Definition
\ref{def.rowscols}, in Definition \ref{def.submatrix.minor}, in Definition
\ref{def.block2x2}, in Definition \ref{def.unrows} and in Definition
\ref{def.sol.det.rk1upd.scalar}.

We begin with some simple lemmas.

\begin{lemma}
\label{lem.sol.det.bordered.1}Let $n\in\mathbb{N}$. Let $u\in\mathbb{K}%
^{n\times1}$ be a column vector with $n$ entries, and let $v=\left(
v_{1},v_{2},\ldots,v_{n}\right)  \in\mathbb{K}^{1\times n}$ be a row vector
with $n$ entries. (Thus, $uv$ is an $n\times n$-matrix, whereas $vu$ is a
$1\times1$-matrix.) Let $h\in\mathbb{K}$. Let $H$ be the $1\times1$-matrix
$\left(
\begin{array}
[c]{c}%
h
\end{array}
\right)  \in\mathbb{K}^{1\times1}$. Let $A\in\mathbb{K}^{n\times n}$ be an
$n\times n$-matrix. Let $C$ be the $\left(  n+1\right)  \times\left(
n+1\right)  $-matrix $\left(
\begin{array}
[c]{cc}%
A & u\\
v & H
\end{array}
\right)  $. Write this $\left(  n+1\right)  \times\left(  n+1\right)  $-matrix
$C$ in the form $C=\left(  c_{i,j}\right)  _{1\leq i\leq n+1,\ 1\leq j\leq
n+1}$. Then:

\textbf{(a)} We have $c_{n+1,q}=v_{q}$ for each $q\in\left\{  1,2,\ldots
,n\right\}  $.

\textbf{(b)} We have $c_{n+1,n+1}=h$.

\textbf{(c)} We have $C_{\sim\left(  n+1\right)  ,\bullet}=\left(  A\mid
u\right)  $.

\textbf{(d)} We have $C_{\sim\left(  n+1\right)  ,\sim q}=\left(  A\mid
u\right)  _{\bullet,\sim q}$ for each $q\in\left\{  1,2,\ldots,n+1\right\}  $.

\textbf{(e)} We have $C_{\sim\left(  n+1\right)  ,\sim\left(  n+1\right)  }=A$.
\end{lemma}

\begin{vershort}
\begin{proof}
[Proof of Lemma \ref{lem.sol.det.bordered.1}.]All parts of Lemma
\ref{lem.sol.det.bordered.1} are obvious from a look at the matrix $C=\left(
\begin{array}
[c]{cc}%
A & u\\
v & H
\end{array}
\right)  $: For example, part \textbf{(d)} says that removing the $\left(
n+1\right)  $-st row and the $q$-th column (for some $q\in\left\{
1,2,\ldots,n+1\right\}  $) from this matrix yields the same result as
attaching the column vector $u$ to $A$ at its right and then removing the
$q$-th column (which should be obvious). Turning this into a rigorous proof is
a straightforward exercise in bookkeeping (using (\ref{eq.def.block2x2.formal}%
), Proposition \ref{prop.unrows.basics} and Proposition
\ref{prop.addcol.props1}).
\end{proof}
\end{vershort}

\begin{verlong}
\begin{proof}
[Proof of Lemma \ref{lem.sol.det.bordered.1}.]The definition of $C$ yields
\begin{equation}
C=\left(
\begin{array}
[c]{cc}%
A & u\\
v & H
\end{array}
\right)  . \label{pf.lem.sol.det.bordered.1.C=}%
\end{equation}

Write the column vector $u\in\mathbb{K}^{n\times1}$ in the form $u=\left(
\begin{array}
[c]{c}%
u_{1}\\
u_{2}\\
\vdots\\
u_{n}%
\end{array}
\right)  $. Thus, $u=\left(
\begin{array}
[c]{c}%
u_{1}\\
u_{2}\\
\vdots\\
u_{n}%
\end{array}
\right)  =\left(  u_{1},u_{2},\ldots,u_{n}\right)  ^{T}=\left(  u_{i}\right)
_{1\leq i\leq n,\ 1\leq j\leq1}$.

Also, $v=\left(  v_{1},v_{2},\ldots,v_{n}\right)  =\left(  v_{j}\right)
_{1\leq i\leq1,\ 1\leq j\leq n}$ and $H=\left(
\begin{array}
[c]{c}%
h
\end{array}
\right)  =\left(  h\right)  _{1\leq i\leq1,\ 1\leq j\leq1}$.

Write the matrix $A\in\mathbb{K}^{n\times n}$ in the form $A=\left(
a_{i,j}\right)  _{1\leq i\leq n,\ 1\leq j\leq n}$.

Now, we have $A=\left(  a_{i,j}\right)  _{1\leq i\leq n,\ 1\leq j\leq n}$,
$u=\left(  u_{i}\right)  _{1\leq i\leq n,\ 1\leq j\leq1}$, $v=\left(
v_{j}\right)  _{1\leq i\leq1,\ 1\leq j\leq n}$ and $H=\left(  h\right)
_{1\leq i\leq1,\ 1\leq j\leq1}$. Thus, (\ref{eq.def.block2x2.formal}) (applied
to $n$, $1$, $n$, $1$, $A$, $u$, $v$, $H$, $a_{i,j}$, $u_{i}$, $v_{j}$ and $h$
instead of $n$, $n^{\prime}$, $m$, $m^{\prime}$, $A$, $B$, $C$, $D$, $a_{i,j}%
$, $b_{i,j}$, $c_{i,j}$ and $d_{i,j}$) yields%
\[
\left(
\begin{array}
[c]{cc}%
A & u\\
v & H
\end{array}
\right)  =\left(
\begin{cases}
a_{i,j}, & \text{if }i\leq n\text{ and }j\leq n;\\
u_{i}, & \text{if }i\leq n\text{ and }j>n;\\
v_{j}, & \text{if }i>n\text{ and }j\leq n;\\
h, & \text{if }i>n\text{ and }j>n
\end{cases}
\right)  _{1\leq i\leq n+1,\ 1\leq j\leq n+1}.
\]
Comparing this with%
\begin{align*}
\left(
\begin{array}
[c]{cc}%
A & u\\
v & H
\end{array}
\right)   &  =C\ \ \ \ \ \ \ \ \ \ \left(  \text{by
(\ref{pf.lem.sol.det.bordered.1.C=})}\right) \\
&  =\left(  c_{i,j}\right)  _{1\leq i\leq n+1,\ 1\leq j\leq n+1},
\end{align*}
we obtain%
\[
\left(  c_{i,j}\right)  _{1\leq i\leq n+1,\ 1\leq j\leq n+1}=\left(
\begin{cases}
a_{i,j}, & \text{if }i\leq n\text{ and }j\leq n;\\
u_{i}, & \text{if }i\leq n\text{ and }j>n;\\
v_{j}, & \text{if }i>n\text{ and }j\leq n;\\
h, & \text{if }i>n\text{ and }j>n
\end{cases}
\right)  _{1\leq i\leq n+1,\ 1\leq j\leq n+1}.
\]
In other words,%
\begin{equation}
c_{i,j}=%
\begin{cases}
a_{i,j}, & \text{if }i\leq n\text{ and }j\leq n;\\
u_{i}, & \text{if }i\leq n\text{ and }j>n;\\
v_{j}, & \text{if }i>n\text{ and }j\leq n;\\
h, & \text{if }i>n\text{ and }j>n
\end{cases}
\label{pf.lem.sol.det.bordered.1.cij=}%
\end{equation}
for every $\left(  i,j\right)  \in\left\{  1,2,\ldots,n+1\right\}  ^{2}$.

We have $n\in\mathbb{N}$; thus, $n+1$ is a positive integer. Hence,
$n+1\in\left\{  1,2,\ldots,n+1\right\}  $.

\textbf{(a)} Let $q\in\left\{  1,2,\ldots,n\right\}  $. Thus, $q\in\left\{
1,2,\ldots,n\right\}  \subseteq\left\{  1,2,\ldots,n+1\right\}  $. Also, from
$q\in\left\{  1,2,\ldots,n\right\}  $, we obtain $q\leq n$.

Combining $n+1\in\left\{  1,2,\ldots,n+1\right\}  $ with $q\in\left\{
1,2,\ldots,n+1\right\}  $, we obtain $\left(  n+1,q\right)  \in\left\{
1,2,\ldots,n+1\right\}  ^{2}$. Hence, (\ref{pf.lem.sol.det.bordered.1.cij=})
(applied to $\left(  i,j\right)  =\left(  n+1,q\right)  $) yields%
\[
c_{n+1,q}=%
\begin{cases}
a_{n+1,q}, & \text{if }n+1\leq n\text{ and }q\leq n;\\
u_{n+1}, & \text{if }n+1\leq n\text{ and }q>n;\\
v_{q}, & \text{if }n+1>n\text{ and }q\leq n;\\
h, & \text{if }n+1>n\text{ and }q>n
\end{cases}
=v_{q}%
\]
(since $n+1>n$ and $q\leq n$). This proves Lemma \ref{lem.sol.det.bordered.1}
\textbf{(a)}.

\textbf{(b)} Combining $n+1\in\left\{  1,2,\ldots,n+1\right\}  $ with
$n+1\in\left\{  1,2,\ldots,n+1\right\}  $, we obtain $\left(  n+1,n+1\right)
\in\left\{  1,2,\ldots,n+1\right\}  ^{2}$. Hence,
(\ref{pf.lem.sol.det.bordered.1.cij=}) (applied to $\left(  i,j\right)
=\left(  n+1,n+1\right)  $) yields%
\[
c_{n+1,n+1}=%
\begin{cases}
a_{n+1,n+1}, & \text{if }n+1\leq n\text{ and }n+1\leq n;\\
u_{n+1}, & \text{if }n+1\leq n\text{ and }n+1>n;\\
v_{n+1}, & \text{if }n+1>n\text{ and }n+1\leq n;\\
h, & \text{if }n+1>n\text{ and }n+1>n
\end{cases}
=h
\]
(since $n+1>n$ and $n+1>n$). This proves Lemma \ref{lem.sol.det.bordered.1}
\textbf{(b)}.

\textbf{(c)} The definition of $C_{\sim\left(  n+1\right)  ,\bullet}$ yields%
\begin{align*}
C_{\sim\left(  n+1\right)  ,\bullet}  &  =\operatorname*{rows}%
\nolimits_{1,2,\ldots,\widehat{n+1},\ldots,n+1}C\ \ \ \ \ \ \ \ \ \ \left(
\text{since }C\text{ is an }\left(  n+1\right)  \times\left(  n+1\right)
\text{-matrix}\right) \\
&  =\operatorname*{rows}\nolimits_{1,2,\ldots,n}C\ \ \ \ \ \ \ \ \ \ \left(
\text{since }\left(  1,2,\ldots,\widehat{n+1},\ldots,n+1\right)  =\left(
1,2,\ldots,n\right)  \right) \\
&  =\left(  c_{x,j}\right)  _{1\leq x\leq n,\ 1\leq j\leq n+1}%
\ \ \ \ \ \ \ \ \ \ \left(
\begin{array}
[c]{c}%
\text{by the definition of }\operatorname*{rows}\nolimits_{1,2,\ldots,n}C\\
\text{(since }C=\left(  c_{i,j}\right)  _{1\leq i\leq n+1,\ 1\leq j\leq
n+1}\text{)}%
\end{array}
\right) \\
&  =\left(  c_{i,j}\right)  _{1\leq i\leq n,\ 1\leq j\leq n+1}%
\end{align*}
(here, we have renamed the index $\left(  x,j\right)  $ as $\left(
i,j\right)  $). In particular, $C_{\sim\left(  n+1\right)  ,\bullet}%
\in\mathbb{K}^{n\times\left(  n+1\right)  }$.

We shall now prove that
\begin{equation}
\left(  C_{\sim\left(  n+1\right)  ,\bullet}\right)  _{\bullet,q}=\left(
A\mid u\right)  _{\bullet,q}\ \ \ \ \ \ \ \ \ \ \text{for each }q\in\left\{
1,2,\ldots,n+1\right\}  . \label{pf.lem.sol.det.bordered.1.c.goal}%
\end{equation}

[\textit{Proof of (\ref{pf.lem.sol.det.bordered.1.c.goal}):} Let $q\in\left\{
1,2,\ldots,n+1\right\}  $.

Proposition \ref{prop.unrows.basics} \textbf{(b)} (applied to $n+1$,
$C_{\sim\left(  n+1\right)  ,\bullet}$ and $q$ instead of $m$, $A$ and $v$)
yields%
\[
\left(  C_{\sim\left(  n+1\right)  ,\bullet}\right)  _{\bullet,q}=\left(
\text{the }q\text{-th column of the matrix }C_{\sim\left(  n+1\right)
,\bullet}\right)  =\operatorname*{cols}\nolimits_{q}\left(  C_{\sim\left(
n+1\right)  ,\bullet}\right)  .
\]
Hence,%
\begin{align}
\left(  C_{\sim\left(  n+1\right)  ,\bullet}\right)  _{\bullet,q}  &  =\left(
\text{the }q\text{-th column of the matrix }\underbrace{C_{\sim\left(
n+1\right)  ,\bullet}}_{=\left(  c_{i,j}\right)  _{1\leq i\leq n,\ 1\leq j\leq
n+1}}\right) \nonumber\\
&  =\left(  \text{the }q\text{-th column of the matrix }\left(  c_{i,j}%
\right)  _{1\leq i\leq n,\ 1\leq j\leq n+1}\right) \nonumber\\
&  =\left(  c_{i,q}\right)  _{1\leq i\leq n,\ 1\leq j\leq1}.
\label{pf.lem.sol.det.bordered.1.c.goal.pf.1}%
\end{align}

But we are in one of the following two cases:

\textit{Case 1:} We have $q\neq n+1$.

\textit{Case 2:} We have $q=n+1$.

Let us first consider Case 1. In this case, we have $q\neq n+1$. Combining
$q\in\left\{  1,2,\ldots,n+1\right\}  $ with $q\neq n+1$, we obtain
$q\in\left\{  1,2,\ldots,n+1\right\}  \setminus\left\{  n+1\right\}  =\left\{
1,2,\ldots,n\right\}  $. Hence, Proposition \ref{prop.addcol.props1}
\textbf{(a)} (applied to $n$ and $u$ instead of $m$ and $v$) yields
\begin{align}
\left(  A\mid u\right)  _{\bullet,q}  &  =A_{\bullet,q}=\left(  \text{the
}q\text{-th column of the matrix }\underbrace{A}_{=\left(  a_{i,j}\right)
_{1\leq i\leq n,\ 1\leq j\leq n}}\right) \nonumber\\
&  =\left(  \text{the }q\text{-th column of the matrix }\left(  a_{i,j}%
\right)  _{1\leq i\leq n,\ 1\leq j\leq n}\right) \nonumber\\
&  =\left(  a_{i,q}\right)  _{1\leq i\leq n,\ 1\leq j\leq1}.
\label{pf.lem.sol.det.bordered.1.c.goal.pf.c1.2}%
\end{align}

But each $i\in\left\{  1,2,\ldots,n\right\}  $ satisfies $c_{i,q}=a_{i,q}%
$\ \ \ \ \footnote{\textit{Proof.} Let $i\in\left\{  1,2,\ldots,n\right\}  $.
Thus, $i\leq n$. Also, $q\leq n$ (since $q\in\left\{  1,2,\ldots,n\right\}
$). From $i\in\left\{  1,2,\ldots,n\right\}  \subseteq\left\{  1,2,\ldots
,n+1\right\}  $ and $q\in\left\{  1,2,\ldots,n+1\right\}  $, we obtain
$\left(  i,q\right)  \in\left\{  1,2,\ldots,n+1\right\}  ^{2}$. Thus,
(\ref{pf.lem.sol.det.bordered.1.cij=}) (applied to $\left(  i,q\right)  $
instead of $\left(  i,j\right)  $) yields%
\[
c_{i,q}=%
\begin{cases}
a_{i,q}, & \text{if }i\leq n\text{ and }q\leq n;\\
u_{i}, & \text{if }i\leq n\text{ and }q>n;\\
v_{q}, & \text{if }i>n\text{ and }q\leq n;\\
h, & \text{if }i>n\text{ and }q>n
\end{cases}
=a_{i,q}%
\]
(since $i\leq n$ and $q\leq n$). Qed.}. Hence, $\left(  \underbrace{c_{i,q}%
}_{=a_{i,q}}\right)  _{1\leq i\leq n,\ 1\leq j\leq1}=\left(  a_{i,q}\right)
_{1\leq i\leq n,\ 1\leq j\leq1}$. Thus,
(\ref{pf.lem.sol.det.bordered.1.c.goal.pf.1}) becomes%
\[
\left(  C_{\sim\left(  n+1\right)  ,\bullet}\right)  _{\bullet,q}=\left(
c_{i,q}\right)  _{1\leq i\leq n,\ 1\leq j\leq1}=\left(  a_{i,q}\right)
_{1\leq i\leq n,\ 1\leq j\leq1}.
\]
Comparing this with (\ref{pf.lem.sol.det.bordered.1.c.goal.pf.c1.2}), we
obtain $\left(  C_{\sim\left(  n+1\right)  ,\bullet}\right)  _{\bullet
,q}=\left(  A\mid u\right)  _{\bullet,q}$. Hence,
(\ref{pf.lem.sol.det.bordered.1.c.goal}) is proven in Case 1.

Let us now consider Case 2. In this case, we have $q=n+1$. But Proposition
\ref{prop.addcol.props1} \textbf{(b)} (applied to $n$ and $u$ instead of $m$
and $v$) yields $\left(  A\mid u\right)  _{\bullet,n+1}=u$. From $q=n+1$, we
obtain
\begin{equation}
\left(  A\mid u\right)  _{\bullet,q}=\left(  A\mid u\right)  _{\bullet,n+1}=u.
\label{pf.lem.sol.det.bordered.1.c.goal.pf.c2.2}%
\end{equation}

On the other hand, each $i\in\left\{  1,2,\ldots,n\right\}  $ satisfies
$c_{i,q}=u_{i}$\ \ \ \ \footnote{\textit{Proof.} Let $i\in\left\{
1,2,\ldots,n\right\}  $. Thus, $i\leq n$. Also, $q=n+1>n$. From $i\in\left\{
1,2,\ldots,n\right\}  \subseteq\left\{  1,2,\ldots,n+1\right\}  $ and
$q\in\left\{  1,2,\ldots,n+1\right\}  $, we obtain $\left(  i,q\right)
\in\left\{  1,2,\ldots,n+1\right\}  ^{2}$. Thus,
(\ref{pf.lem.sol.det.bordered.1.cij=}) (applied to $\left(  i,q\right)  $
instead of $\left(  i,j\right)  $) yields%
\[
c_{i,q}=%
\begin{cases}
a_{i,q}, & \text{if }i\leq n\text{ and }q\leq n;\\
u_{i}, & \text{if }i\leq n\text{ and }q>n;\\
v_{q}, & \text{if }i>n\text{ and }q\leq n;\\
h, & \text{if }i>n\text{ and }q>n
\end{cases}
=u_{i}%
\]
(since $i\leq n$ and $q>n$). Qed.}. Hence,
\[
\left(  \underbrace{c_{i,q}}_{=u_{i}}\right)  _{1\leq i\leq n,\ 1\leq j\leq
1}=\left(  u_{i}\right)  _{1\leq i\leq n,\ 1\leq j\leq1}=u
\]
(since $u=\left(  u_{i}\right)  _{1\leq i\leq n,\ 1\leq j\leq1}$). Thus,
(\ref{pf.lem.sol.det.bordered.1.c.goal.pf.1}) becomes%
\[
\left(  C_{\sim\left(  n+1\right)  ,\bullet}\right)  _{\bullet,q}=\left(
c_{i,q}\right)  _{1\leq i\leq n,\ 1\leq j\leq1}=u.
\]
Comparing this with (\ref{pf.lem.sol.det.bordered.1.c.goal.pf.c2.2}), we
obtain $\left(  C_{\sim\left(  n+1\right)  ,\bullet}\right)  _{\bullet
,q}=\left(  A\mid u\right)  _{\bullet,q}$. Hence,
(\ref{pf.lem.sol.det.bordered.1.c.goal}) is proven in Case 2.

We have now proven (\ref{pf.lem.sol.det.bordered.1.c.goal}) in each of the two
Cases 1 and 2. Since these two Cases cover all possibilities, we thus conclude
that (\ref{pf.lem.sol.det.bordered.1.c.goal}) always holds.]

We have now proven that $\left(  C_{\sim\left(  n+1\right)  ,\bullet}\right)
_{\bullet,q}=\left(  A\mid u\right)  _{\bullet,q}$ for each $q\in\left\{
1,2,\ldots,n+1\right\}  $. Hence, Lemma \ref{lem.sol.prop.addcol.props.cols}
(applied to $n+1$, $C_{\sim\left(  n+1\right)  ,\bullet}$ and $\left(  A\mid
u\right)  $ instead of $m$, $A$ and $B$) shows that $C_{\sim\left(
n+1\right)  ,\bullet}=\left(  A\mid u\right)  $ (since $C_{\sim\left(
n+1\right)  ,\bullet}$ and $\left(  A\mid u\right)  $ are two $n\times\left(
n+1\right)  $-matrices). This proves Lemma \ref{lem.sol.det.bordered.1}
\textbf{(c)}.

\textbf{(d)} Let $q\in\left\{  1,2,\ldots,n+1\right\}  $. Recall that
$n+1\in\left\{  1,2,\ldots,n+1\right\}  $.

Thus, Proposition \ref{prop.unrows.basics} \textbf{(c)} (applied to $n+1$,
$n+1$, $C$, $n+1$ and $q$ instead of $n$, $m$, $A$, $u$ and $v$) shows that%
\[
\left(  C_{\bullet,\sim q}\right)  _{\sim\left(  n+1\right)  ,\bullet}=\left(
C_{\sim\left(  n+1\right)  ,\bullet}\right)  _{\bullet,\sim q}=C_{\sim\left(
n+1\right)  ,\sim q}.
\]
Hence,
\[
C_{\sim\left(  n+1\right)  ,\sim q}=\left(  \underbrace{C_{\sim\left(
n+1\right)  ,\bullet}}_{\substack{=\left(  A\mid u\right)  \\\text{(by Lemma
\ref{lem.sol.det.bordered.1} \textbf{(c)})}}}\right)  _{\bullet,\sim
q}=\left(  A\mid u\right)  _{\bullet,\sim q}.
\]
This proves Lemma \ref{lem.sol.det.bordered.1} \textbf{(d)}.

\textbf{(e)} We have $n+1\in\left\{  1,2,\ldots,n+1\right\}  $. Thus, Lemma
\ref{lem.sol.det.bordered.1} \textbf{(d)} (applied to $q=n+1$) yields
$C_{\sim\left(  n+1\right)  ,\sim\left(  n+1\right)  }=\left(  A\mid u\right)
_{\bullet,\sim\left(  n+1\right)  }=A$ (by Proposition
\ref{prop.addcol.props1} \textbf{(d)} (applied to $n$ and $u$ instead of $m$
and $v$)). This proves Lemma \ref{lem.sol.det.bordered.1} \textbf{(e)}.
\end{proof}
\end{verlong}

\begin{lemma}
\label{lem.sol.det.bordered.2}Let $n\in\mathbb{N}$. Let $u=\left(  u_{1}%
,u_{2},\ldots,u_{n}\right)  ^{T}\in\mathbb{K}^{n\times1}$ be a column vector
with $n$ entries, and let $v=\left(  v_{1},v_{2},\ldots,v_{n}\right)
\in\mathbb{K}^{1\times n}$ be a row vector with $n$ entries. Let
$h\in\mathbb{K}$. Let $H$ be the $1\times1$-matrix $\left(
\begin{array}
[c]{c}%
h
\end{array}
\right)  \in\mathbb{K}^{1\times1}$. Let $A\in\mathbb{K}^{n\times n}$ be an
$n\times n$-matrix. Then,%
\[
\det\left(
\begin{array}
[c]{cc}%
A & u\\
v & H
\end{array}
\right)  =h\det A-\sum_{i=1}^{n}\sum_{j=1}^{n}\left(  -1\right)  ^{i+j}%
u_{j}v_{i}\det\left(  A_{\sim j,\sim i}\right)  .
\]

\end{lemma}

\begin{proof}
[Proof of Lemma \ref{lem.sol.det.bordered.2}.]Let $C$ be the $\left(
n+1\right)  \times\left(  n+1\right)  $-matrix $\left(
\begin{array}
[c]{cc}%
A & u\\
v & H
\end{array}
\right)  $. Write this $\left(  n+1\right)  \times\left(  n+1\right)  $-matrix
$C$ in the form $C=\left(  c_{i,j}\right)  _{1\leq i\leq n+1,\ 1\leq j\leq
n+1}$.

\begin{verlong}
We have $n\in\mathbb{N}$; thus, $n+1$ is a positive integer. Hence,
$n+1\in\left\{  1,2,\ldots,n+1\right\}  $.
\end{verlong}

Thus, Theorem \ref{thm.laplace.gen} \textbf{(a)} (applied to $n+1$, $C$,
$c_{i,j}$ and $n+1$ instead of $n$, $A$, $a_{i,j}$ and $n+1$) yields%
\begin{align*}
&  \det C\\
&  =\sum_{q=1}^{n+1}\left(  -1\right)  ^{\left(  n+1\right)  +q}c_{n+1,q}%
\det\left(  C_{\sim\left(  n+1\right)  ,\sim q}\right) \\
&  \ \ \ \ \ \ \ \ \ \ \left(  \text{since }C=\left(  c_{i,j}\right)  _{1\leq
i\leq n+1,\ 1\leq j\leq n+1}\right) \\
&  =\sum_{q=1}^{n}\left(  -1\right)  ^{\left(  n+1\right)  +q}%
\underbrace{c_{n+1,q}}_{\substack{=v_{q}\\\text{(by Lemma
\ref{lem.sol.det.bordered.1} \textbf{(a)})}}}\det\left(  \underbrace{C_{\sim
\left(  n+1\right)  ,\sim q}}_{\substack{=\left(  A\mid u\right)
_{\bullet,\sim q}\\\text{(by Lemma \ref{lem.sol.det.bordered.1} \textbf{(d)}%
}\\\text{(since }q\in\left\{  1,2,\ldots,n\right\}  \subseteq\left\{
1,2,\ldots,n+1\right\}  \text{))}}}\right) \\
&  \ \ \ \ \ \ \ \ \ \ +\underbrace{\left(  -1\right)  ^{\left(  n+1\right)
+\left(  n+1\right)  }}_{\substack{=1\\\text{(since }\left(  n+1\right)
+\left(  n+1\right)  =2\left(  n+1\right)  \text{ is even)}}%
}\ \ \underbrace{c_{n+1,n+1}}_{\substack{=h\\\text{(by Lemma
\ref{lem.sol.det.bordered.1} \textbf{(b)})}}}\det\left(  \underbrace{C_{\sim
\left(  n+1\right)  ,\sim\left(  n+1\right)  }}_{\substack{=A\\\text{(by Lemma
\ref{lem.sol.det.bordered.1} \textbf{(e)})}}}\right) \\
&  \ \ \ \ \ \ \ \ \ \ \left(  \text{here, we have split off the addend for
}q=n+1\text{ from the sum}\right) \\
&  =\sum_{q=1}^{n}\left(  -1\right)  ^{\left(  n+1\right)  +q}v_{q}\det\left(
\left(  A\mid u\right)  _{\bullet,\sim q}\right)  +h\det A.
\end{align*}

\begin{vershort}
Subtracting $h\det A$ from both sides of this equality, we obtain%
\begin{align*}
&  \det C-h\det A\\
&  =\sum_{q=1}^{n}\left(  -1\right)  ^{\left(  n+1\right)  +q}v_{q}%
\underbrace{\det\left(  \left(  A\mid u\right)  _{\bullet,\sim q}\right)
}_{\substack{=\sum_{j=1}^{n}\left(  -1\right)  ^{n+j}u_{j}\det\left(  A_{\sim
j,\sim q}\right)  \\\text{(by Lemma \ref{lem.sol.det.rk1upd.detAuk} (applied
to }k=q\text{))}}}\\
&  =\sum_{q=1}^{n}\left(  -1\right)  ^{\left(  n+1\right)  +q}v_{q}\left(
\sum_{j=1}^{n}\left(  -1\right)  ^{n+j}u_{j}\det\left(  A_{\sim j,\sim
q}\right)  \right) \\
&  =\sum_{q=1}^{n}\sum_{j=1}^{n}\left(  -1\right)  ^{\left(  n+1\right)
+q}v_{q}\left(  -1\right)  ^{n+j}u_{j}\det\left(  A_{\sim j,\sim q}\right) \\
&  =\sum_{q=1}^{n}\sum_{j=1}^{n}\underbrace{\left(  -1\right)  ^{\left(
n+1\right)  +q}\left(  -1\right)  ^{n+j}}_{\substack{=\left(  -1\right)
^{\left(  \left(  n+1\right)  +q\right)  +\left(  n+j\right)  }=\left(
-1\right)  ^{q+j+1}\\\text{(since }\left(  \left(  n+1\right)  +q\right)
+\left(  n+j\right)  =\left(  q+j+1\right)  +2n\equiv q+j+1\operatorname{mod}%
2\text{)}}}u_{j}v_{q}\det\left(  A_{\sim j,\sim q}\right) \\
&  =\sum_{q=1}^{n}\sum_{j=1}^{n}\underbrace{\left(  -1\right)  ^{q+j+1}%
}_{=-\left(  -1\right)  ^{q+j}}u_{j}v_{q}\det\left(  A_{\sim j,\sim q}\right)
\\
&  =-\sum_{q=1}^{n}\sum_{j=1}^{n}\left(  -1\right)  ^{q+j}u_{j}v_{q}%
\det\left(  A_{\sim j,\sim q}\right)  =-\sum_{i=1}^{n}\sum_{j=1}^{n}\left(
-1\right)  ^{i+j}u_{j}v_{i}\det\left(  A_{\sim j,\sim i}\right)
\end{align*}
(here, we have renamed the summation index $q$ as $i$).
\end{vershort}

\begin{verlong}
Subtracting $h\det A$ from both sides of this equality, we obtain%
\begin{align*}
&  \det C-h\det A\\
&  =\sum_{q=1}^{n}\left(  -1\right)  ^{\left(  n+1\right)  +q}v_{q}%
\underbrace{\det\left(  \left(  A\mid u\right)  _{\bullet,\sim q}\right)
}_{\substack{=\sum_{j=1}^{n}\left(  -1\right)  ^{n+j}u_{j}\det\left(  A_{\sim
j,\sim q}\right)  \\\text{(by Lemma \ref{lem.sol.det.rk1upd.detAuk} (applied
to }k=q\text{))}}}\\
&  =\sum_{q=1}^{n}\left(  -1\right)  ^{\left(  n+1\right)  +q}v_{q}\left(
\sum_{j=1}^{n}\left(  -1\right)  ^{n+j}u_{j}\det\left(  A_{\sim j,\sim
q}\right)  \right) \\
&  =\sum_{q=1}^{n}\sum_{j=1}^{n}\left(  -1\right)  ^{\left(  n+1\right)
+q}\underbrace{v_{q}\left(  -1\right)  ^{n+j}u_{j}}_{=\left(  -1\right)
^{n+j}u_{j}v_{q}}\det\left(  A_{\sim j,\sim q}\right) \\
&  =\sum_{q=1}^{n}\sum_{j=1}^{n}\underbrace{\left(  -1\right)  ^{\left(
n+1\right)  +q}\left(  -1\right)  ^{n+j}}_{\substack{=\left(  -1\right)
^{\left(  \left(  n+1\right)  +q\right)  +\left(  n+j\right)  }=\left(
-1\right)  ^{q+j+1}\\\text{(since }\left(  \left(  n+1\right)  +q\right)
+\left(  n+j\right)  =\left(  q+j+1\right)  +2n\equiv q+j+1\operatorname{mod}%
2\text{)}}}u_{j}v_{q}\det\left(  A_{\sim j,\sim q}\right) \\
&  =\sum_{q=1}^{n}\sum_{j=1}^{n}\underbrace{\left(  -1\right)  ^{q+j+1}%
}_{=-\left(  -1\right)  ^{q+j}}u_{j}v_{q}\det\left(  A_{\sim j,\sim q}\right)
\\
&  =\sum_{q=1}^{n}\sum_{j=1}^{n}\left(  -\left(  -1\right)  ^{q+j}\right)
u_{j}v_{q}\det\left(  A_{\sim j,\sim q}\right)  =-\sum_{q=1}^{n}\sum_{j=1}%
^{n}\left(  -1\right)  ^{q+j}u_{j}v_{q}\det\left(  A_{\sim j,\sim q}\right) \\
&  =-\sum_{i=1}^{n}\sum_{j=1}^{n}\left(  -1\right)  ^{i+j}u_{j}v_{i}%
\det\left(  A_{\sim j,\sim i}\right)
\end{align*}
(here, we have renamed the summation index $q$ as $i$).
\end{verlong}

Adding $h\det A$ to both sides of this equality, we obtain%
\[
\det C=h\det A-\sum_{i=1}^{n}\sum_{j=1}^{n}\left(  -1\right)  ^{i+j}u_{j}%
v_{i}\det\left(  A_{\sim j,\sim i}\right)  .
\]
In view of $C=\left(
\begin{array}
[c]{cc}%
A & u\\
v & H
\end{array}
\right)  $ (by the definition of $C$), this rewrites as%
\[
\det\left(
\begin{array}
[c]{cc}%
A & u\\
v & H
\end{array}
\right)  =h\det A-\sum_{i=1}^{n}\sum_{j=1}^{n}\left(  -1\right)  ^{i+j}%
u_{j}v_{i}\det\left(  A_{\sim j,\sim i}\right)  .
\]
This proves Lemma \ref{lem.sol.det.bordered.2}.
\end{proof}

We can now easily obtain the following corollary, which is just a restatement
of Exercise \ref{exe.det.bordered} \textbf{(a)}:

\begin{corollary}
\label{cor.sol.det.bordered.3}Let $n\in\mathbb{N}$. Let $u\in\mathbb{K}%
^{n\times1}$ be a column vector with $n$ entries, and let $v\in\mathbb{K}%
^{1\times n}$ be a row vector with $n$ entries. (Thus, $uv$ is an $n\times
n$-matrix, whereas $vu$ is a $1\times1$-matrix.) Let $h\in\mathbb{K}$. Let $H$
be the $1\times1$-matrix $\left(
\begin{array}
[c]{c}%
h
\end{array}
\right)  \in\mathbb{K}^{1\times1}$. Let $A\in\mathbb{K}^{n\times n}$ be an
$n\times n$-matrix. Then,%
\[
\det\left(
\begin{array}
[c]{cc}%
A & u\\
v & H
\end{array}
\right)  =h\det A-\operatorname*{ent}\left(  v\left(  \operatorname*{adj}%
A\right)  u\right)  .
\]
(Here, we are using the notation from Definition
\ref{def.sol.det.rk1upd.scalar}.)
\end{corollary}

\begin{proof}
[Proof of Corollary \ref{cor.sol.det.bordered.3}.]Write the column vector
$u\in\mathbb{K}^{n\times1}$ in the form $u=\left(
\begin{array}
[c]{c}%
u_{1}\\
u_{2}\\
\vdots\\
u_{n}%
\end{array}
\right)  $. Thus, $u=\left(
\begin{array}
[c]{c}%
u_{1}\\
u_{2}\\
\vdots\\
u_{n}%
\end{array}
\right)  =\left(  u_{1},u_{2},\ldots,u_{n}\right)  ^{T}$.

Write the row vector $v\in\mathbb{K}^{1\times n}$ in the form $v=\left(
v_{1},v_{2},\ldots,v_{n}\right)  $. Thus, Corollary
\ref{cor.sol.det.rk1upd.claim.ent(v(adjA)u)} yields%
\begin{equation}
\operatorname*{ent}\left(  v\left(  \operatorname*{adj}A\right)  u\right)
=\sum_{i=1}^{n}\sum_{j=1}^{n}\left(  -1\right)  ^{i+j}u_{j}v_{i}\det\left(
A_{\sim j,\sim i}\right)  . \label{pf.cor.sol.det.bordered.3.1}%
\end{equation}

But Lemma \ref{lem.sol.det.bordered.2} yields%
\begin{align*}
\det\left(
\begin{array}
[c]{cc}%
A & u\\
v & H
\end{array}
\right)   &  =h\det A-\underbrace{\sum_{i=1}^{n}\sum_{j=1}^{n}\left(
-1\right)  ^{i+j}u_{j}v_{i}\det\left(  A_{\sim j,\sim i}\right)
}_{\substack{=\operatorname*{ent}\left(  v\left(  \operatorname*{adj}A\right)
u\right)  \\\text{(by (\ref{pf.cor.sol.det.bordered.3.1}))}}}\\
&  =h\det A-\operatorname*{ent}\left(  v\left(  \operatorname*{adj}A\right)
u\right)  .
\end{align*}
This proves Corollary \ref{cor.sol.det.bordered.3}.
\end{proof}

We now take aim at part \textbf{(b)} of Exercise \ref{exe.det.bordered}. We
begin with a lemma that computes the $\left(  n-1\right)  \times\left(
n-1\right)  $-minors of a diagonal matrix:

\begin{lemma}
\label{lem.det.bordered.det-adj}Let $n\in\mathbb{N}$. For every two objects
$i$ and $j$, define $\delta_{i,j}\in\mathbb{K}$ by $\delta_{i,j}=%
\begin{cases}
1, & \text{if }i=j;\\
0, & \text{if }i\neq j
\end{cases}
$.

Let $d_{1},d_{2},\ldots,d_{n}$ be $n$ elements of $\mathbb{K}$. Let $D$ be the
$n\times n$-matrix $\left(  d_{i}\delta_{i,j}\right)  _{1\leq i\leq n,\ 1\leq
j\leq n}$.

Let $p\in\left\{  1,2,\ldots,n\right\}  $ and $q\in\left\{  1,2,\ldots
,n\right\}  $. Then,
\[
\det\left(  D_{\sim p,\sim q}\right)  =\delta_{p,q}\prod_{\substack{j\in
\left\{  1,2,\ldots,n\right\}  ;\\j\neq p}}d_{j}.
\]

\end{lemma}

\begin{proof}
[Proof of Lemma \ref{lem.det.bordered.det-adj}.]We shall use the notations
from Definition \ref{def.sect.laplace.notations}.

Define a subset $P$ of $\left\{  1,2,\ldots,n\right\}  $ by $P=\left\{
1,2,\ldots,n\right\}  \setminus\left\{  p\right\}  $.

Define a subset $Q$ of $\left\{  1,2,\ldots,n\right\}  $ by $Q=\left\{
1,2,\ldots,n\right\}  \setminus\left\{  q\right\}  $.

\begin{vershort}
Clearly, $P=Q$ holds if and only if $p=q$ holds. Thus, $\delta_{P,Q}%
=\delta_{p,q}$. Also, clearly, $\left\vert P\right\vert =\left\vert
Q\right\vert =n-1$.
\end{vershort}

\begin{verlong}
It is easy to see that $\delta_{P,Q}=\delta_{p,q}$%
\ \ \ \ \footnote{\textit{Proof.} We are in one of the following two cases:
\par
\textit{Case 1:} We have $p\neq q$.
\par
\textit{Case 2:} We have $p=q$.
\par
Let us first consider Case 1. In this case, we have $p\neq q$. Hence,
$\delta_{p,q}=0$. From $p\neq q$, we obtain $p\notin\left\{  q\right\}  $.
Combining $p\in\left\{  1,2,\ldots,n\right\}  $ with $p\notin\left\{
q\right\}  $, we find $p\in\left\{  1,2,\ldots,n\right\}  \setminus\left\{
q\right\}  =Q$ (since $Q=\left\{  1,2,\ldots,n\right\}  \setminus\left\{
q\right\}  $). Since $\left\{  p\right\}  $ is a subset of $\left\{
1,2,\ldots,n\right\}  $ (because of $p\in\left\{  1,2,\ldots,n\right\}  $), we
have
\[
\left\{  1,2,\ldots,n\right\}  \setminus\left(  \left\{  1,2,\ldots,n\right\}
\setminus\left\{  p\right\}  \right)  =\left\{  p\right\}  .
\]
Thus,%
\[
\left\{  p\right\}  =\left\{  1,2,\ldots,n\right\}  \setminus
\underbrace{\left(  \left\{  1,2,\ldots,n\right\}  \setminus\left\{
p\right\}  \right)  }_{=P}=\left\{  1,2,\ldots,n\right\}  \setminus P.
\]
Hence, $p\in\left\{  p\right\}  =\left\{  1,2,\ldots,n\right\}  \setminus P$.
In other words, $p\in\left\{  1,2,\ldots,n\right\}  $ and $p\notin P$. If we
had $Q=P$, then we would have $p\in Q=P$, which would contradict $p\notin P$.
Hence, we cannot have $Q=P$. Thus, we have $Q\neq P$. In other words, $P\neq
Q$. Thus, $\delta_{P,Q}=0$. Comparing this with $\delta_{p,q}=0$, we obtain
$\delta_{P,Q}=\delta_{p,q}$. Hence, $\delta_{P,Q}=\delta_{p,q}$ is proven in
Case 1.
\par
Let us now consider Case 2. In this case, $p=q$. Thus, $\delta_{p,q}=1$. But
\[
P=\left\{  1,2,\ldots,n\right\}  \setminus\left\{  \underbrace{p}%
_{=q}\right\}  =\left\{  1,2,\ldots,n\right\}  \setminus\left\{  q\right\}
=Q
\]
(since $Q=\left\{  1,2,\ldots,n\right\}  \setminus\left\{  q\right\}  $).
Hence, $\delta_{P,Q}=1$. Comparing this with $\delta_{p,q}=1$, we obtain
$\delta_{P,Q}=\delta_{p,q}$. Hence, $\delta_{P,Q}=\delta_{p,q}$ is proven in
Case 2.
\par
We now have proven $\delta_{P,Q}=\delta_{p,q}$ in each of the two Cases 1 and
2. Since these two Cases cover all possibilities, we thus conclude that
$\delta_{P,Q}=\delta_{p,q}$ always holds. Qed.}.

Also, $\left\vert P\right\vert =n-1$\ \ \ \ \footnote{\textit{Proof.} We have
$P=\left\{  1,2,\ldots,n\right\}  \setminus\left\{  p\right\}  $. Thus,%
\begin{align*}
\left\vert P\right\vert  &  =\left\vert \left\{  1,2,\ldots,n\right\}
\setminus\left\{  p\right\}  \right\vert =\underbrace{\left\vert \left\{
1,2,\ldots,n\right\}  \right\vert }_{=n}-1\ \ \ \ \ \ \ \ \ \ \left(
\text{since }p\in\left\{  1,2,\ldots,n\right\}  \right) \\
&  =n-1.
\end{align*}
}. The same argument (applied to $q$ and $Q$ instead of $p$ and $P$) yields
$\left\vert Q\right\vert =n-1$. Hence, $\left\vert P\right\vert
=n-1=\left\vert Q\right\vert $.
\end{verlong}

The definition of $D_{\sim p,\sim q}$ yields
\begin{equation}
D_{\sim p,\sim q}=\operatorname*{sub}\nolimits_{1,2,\ldots,\widehat{p}%
,\ldots,n}^{1,2,\ldots,\widehat{q},\ldots,n}D.
\label{pf.lem.det.bordered.det-adj.2}%
\end{equation}

\begin{vershort}
The definition of the list $\left(  1,2,\ldots,\widehat{p},\ldots,n\right)  $
yields
\[
\left(  1,2,\ldots,\widehat{p},\ldots,n\right)  =\left(  1,2,\ldots
,p-1,p+1,p+2,\ldots,n\right)  .
\]
Meanwhile, the definition of $w\left(  P\right)  $ shows that $w\left(
P\right)  $ is the list of all elements of $P$ in increasing order (with no
repetitions); thus, $w\left(  P\right)  =\left(  1,2,\ldots,p-1,p+1,p+2,\ldots
,n\right)  $ (because $P=\left\{  1,2,\ldots,n\right\}  \setminus\left\{
p\right\}  =\left\{  1,2,\ldots,p-1,p+1,p+2,\ldots,n\right\}  $). Comparing
these two equalities, we obtain $w\left(  P\right)  =\left(  1,2,\ldots
,\widehat{p},\ldots,n\right)  $. Similarly, $w\left(  Q\right)  =\left(
1,2,\ldots,\widehat{q},\ldots,n\right)  $. Using the preceding two equalities,
we find%
\[
\operatorname*{sub}\nolimits_{w\left(  P\right)  }^{w\left(  Q\right)
}D=\operatorname*{sub}\nolimits_{\left(  1,2,\ldots,\widehat{p},\ldots
,n\right)  }^{\left(  1,2,\ldots,\widehat{q},\ldots,n\right)  }%
D=\operatorname*{sub}\nolimits_{1,2,\ldots,\widehat{p},\ldots,n}%
^{1,2,\ldots,\widehat{q},\ldots,n}D.
\]

\end{vershort}

\begin{verlong}
The definition of the list $\left(  1,2,\ldots,\widehat{p},\ldots,n\right)  $
yields
\begin{equation}
\left(  1,2,\ldots,\widehat{p},\ldots,n\right)  =\left(  1,2,\ldots
,p-1,p+1,p+2,\ldots,n\right)  . \label{pf.lem.det.bordered.det-adj.4}%
\end{equation}

But%
\[
P=\left\{  1,2,\ldots,n\right\}  \setminus\left\{  p\right\}  =\left\{
1,2,\ldots,p-1,p+1,p+2,\ldots,n\right\}
\]
(since $p\in\left\{  1,2,\ldots,n\right\}  $). Thus, the list $\left(
1,2,\ldots,p-1,p+1,p+2,\ldots,n\right)  $ is a list of all elements of $P$.
Moreover, this list $\left(  1,2,\ldots,p-1,p+1,p+2,\ldots,n\right)  $ is
strictly increasing (since $1<2<\cdots<p-1<p+1<p+2<\cdots<n$), and thus its
entries are distinct. Hence, this list $\left(  1,2,\ldots,p-1,p+1,p+2,\ldots
,n\right)  $ is the list of all elements of $P$ in increasing order (with no
repetitions)\footnote{because this list $\left(  1,2,\ldots,p-1,p+1,p+2,\ldots
,n\right)  $ is a list of all elements of $P$, and is strictly increasing, and
its entries are distinct}. In other words,%
\begin{align}
&  \left(  1,2,\ldots,p-1,p+1,p+2,\ldots,n\right) \nonumber\\
&  =\left(  \text{the list of all elements of }P\text{ in increasing order
(with no repetitions)}\right)  . \label{pf.lem.det.bordered.det-adj.3a}%
\end{align}
On the other hand, $w\left(  P\right)  $ is the list of all elements of $P$ in
increasing order (with no repetitions)\footnote{by the definition of $w\left(
P\right)  $}. Thus,%
\begin{align*}
&  w\left(  P\right) \\
&  =\left(  \text{the list of all elements of }P\text{ in increasing order
(with no repetitions)}\right)  .
\end{align*}
Comparing this with (\ref{pf.lem.det.bordered.det-adj.3a}), we obtain
\[
w\left(  P\right)  =\left(  1,2,\ldots,p-1,p+1,p+2,\ldots,n\right)  .
\]
Comparing this with (\ref{pf.lem.det.bordered.det-adj.4}), we find $w\left(
P\right)  =\left(  1,2,\ldots,\widehat{p},\ldots,n\right)  $. The same
argument (applied to $q$ and $Q$ instead of $p$ and $P$) yields $w\left(
Q\right)  =\left(  1,2,\ldots,\widehat{q},\ldots,n\right)  $. Now,%
\begin{align*}
\operatorname*{sub}\nolimits_{w\left(  P\right)  }^{w\left(  Q\right)  }D  &
=\operatorname*{sub}\nolimits_{\left(  1,2,\ldots,\widehat{p},\ldots,n\right)
}^{\left(  1,2,\ldots,\widehat{q},\ldots,n\right)  }%
D\ \ \ \ \ \ \ \ \ \ \left(
\begin{array}
[c]{c}%
\text{since }w\left(  P\right)  =\left(  1,2,\ldots,\widehat{p},\ldots
,n\right) \\
\text{and }w\left(  Q\right)  =\left(  1,2,\ldots,\widehat{q},\ldots,n\right)
\end{array}
\right) \\
&  =\operatorname*{sub}\nolimits_{1,2,\ldots,\widehat{p},\ldots,n}%
^{1,2,\ldots,\widehat{q},\ldots,n}D.
\end{align*}

\end{verlong}

Comparing this with (\ref{pf.lem.det.bordered.det-adj.2}), we find $D_{\sim
p,\sim q}=\operatorname*{sub}\nolimits_{w\left(  P\right)  }^{w\left(
Q\right)  }D$. Hence,%
\begin{align*}
\det\left(  D_{\sim p,\sim q}\right)   &  =\det\left(  \operatorname*{sub}%
\nolimits_{w\left(  P\right)  }^{w\left(  Q\right)  }D\right)
=\underbrace{\delta_{P,Q}}_{=\delta_{p,q}}\underbrace{\prod_{i\in P}%
}_{\substack{=\prod_{i\in\left\{  1,2,\ldots,n\right\}  \setminus\left\{
p\right\}  }\\\text{(since }P=\left\{  1,2,\ldots,n\right\}  \setminus\left\{
p\right\}  \text{)}}}d_{i}\\
&  \ \ \ \ \ \ \ \ \ \ \left(  \text{by Lemma \ref{lem.diag.minors}}\right) \\
&  =\delta_{p,q}\prod_{i\in\left\{  1,2,\ldots,n\right\}  \setminus\left\{
p\right\}  }d_{i}=\delta_{p,q}\underbrace{\prod_{j\in\left\{  1,2,\ldots
,n\right\}  \setminus\left\{  p\right\}  }}_{=\prod_{\substack{j\in\left\{
1,2,\ldots,n\right\}  ;\\j\neq p}}}d_{j}\\
&  \ \ \ \ \ \ \ \ \ \ \left(
\begin{array}
[c]{c}%
\text{here, we have renamed the index }i\text{ as }j\\
\text{in the product}%
\end{array}
\right) \\
&  =\delta_{p,q}\prod_{\substack{j\in\left\{  1,2,\ldots,n\right\}  ;\\j\neq
p}}d_{j}.
\end{align*}
This proves Lemma \ref{lem.det.bordered.det-adj}.
\end{proof}

Notice that Lemma \ref{lem.det.bordered.det-adj} generalizes Lemma
\ref{lem.sol.det.rk1upd.adjI} (because if we set $d_{i}=1$ in Lemma
\ref{lem.det.bordered.det-adj}, then the resulting matrix $D$ is $I_{n}$).

We are now finally able to solve Exercise \ref{exe.det.bordered}:

\begin{proof}
[Solution to Exercise \ref{exe.det.bordered}.]\textbf{(a)} Let $A\in
\mathbb{K}^{n\times n}$ be an $n\times n$-matrix. Then, Corollary
\ref{cor.sol.det.bordered.3} yields%
\begin{equation}
\det\left(
\begin{array}
[c]{cc}%
A & u\\
v & H
\end{array}
\right)  =h\det A-\operatorname*{ent}\left(  v\left(  \operatorname*{adj}%
A\right)  u\right)  \label{sol.det.bordered.a.1}%
\end{equation}
(where we are using the notation from Definition
\ref{def.sol.det.rk1upd.scalar}). But if we regard the $1\times1$-matrix
$v\left(  \operatorname*{adj}A\right)  u$ as an element of $\mathbb{K}$, then
$\operatorname*{ent}\left(  v\left(  \operatorname*{adj}A\right)  u\right)
=v\left(  \operatorname*{adj}A\right)  u$, and therefore
(\ref{sol.det.bordered.a.1}) simplifies as follows:%
\[
\det\left(
\begin{array}
[c]{cc}%
A & u\\
v & H
\end{array}
\right)  =h\det A-\underbrace{\operatorname*{ent}\left(  v\left(
\operatorname*{adj}A\right)  u\right)  }_{=v\left(  \operatorname*{adj}%
A\right)  u}=h\det A-v\left(  \operatorname*{adj}A\right)  u.
\]
This solves Exercise \ref{exe.det.bordered} \textbf{(a)}.

\textbf{(b)} We have $d_{i}\delta_{i,j}=0$ for every $\left(  i,j\right)
\in\left\{  1,2,\ldots,n\right\}  ^{2}$ satisfying $i<j$%
\ \ \ \ \footnote{\textit{Proof.} Let $\left(  i,j\right)  \in\left\{
1,2,\ldots,n\right\}  ^{2}$ be such that $i<j$. From $i<j$, we obtain $i\neq
j$; thus, $\delta_{i,j}=0$. Hence, $d_{i}\underbrace{\delta_{i,j}}_{=0}=0$,
qed.}. Hence, Exercise \ref{exe.ps4.3} (applied to $D$ and $d_{i}\delta_{i,j}$
instead of $A$ and $a_{i,j}$) yields%
\begin{align*}
\det D  &  =\left(  d_{1}\delta_{1,1}\right)  \left(  d_{2}\delta
_{2,2}\right)  \cdots\left(  d_{n}\delta_{n,n}\right)
\ \ \ \ \ \ \ \ \ \ \left(  \text{since }D=\left(  d_{i}\delta_{i,j}\right)
_{1\leq i\leq n,\ 1\leq j\leq n}\right) \\
&  =\prod_{k=1}^{n}\left(  d_{k}\underbrace{\delta_{k,k}}%
_{\substack{=1\\\text{(since }k=k\text{)}}}\right)  =\prod_{k=1}^{n}%
d_{k}=d_{1}d_{2}\cdots d_{n}.
\end{align*}

Applying Lemma \ref{lem.sol.det.bordered.2} to $A=D$, we obtain%
\begin{equation}
\det\left(
\begin{array}
[c]{cc}%
D & u\\
v & H
\end{array}
\right)  =h\det D-\sum_{i=1}^{n}\sum_{j=1}^{n}\left(  -1\right)  ^{i+j}%
u_{j}v_{i}\det\left(  D_{\sim j,\sim i}\right)  . \label{sol.det.bordered.b.2}%
\end{equation}

But each $i\in\left\{  1,2,\ldots,n\right\}  $ satisfies
\begin{align}
&  \sum_{j=1}^{n}\left(  -1\right)  ^{i+j}u_{j}v_{i}\det\left(  D_{\sim j,\sim
i}\right) \nonumber\\
&  =\underbrace{\sum_{p=1}^{n}}_{=\sum_{p\in\left\{  1,2,\ldots,n\right\}  }%
}\left(  -1\right)  ^{i+p}u_{p}v_{i}\underbrace{\det\left(  D_{\sim p,\sim
i}\right)  }_{\substack{=\delta_{p,i}\prod_{\substack{j\in\left\{
1,2,\ldots,n\right\}  ;\\j\neq p}}d_{j}\\\text{(by Lemma
\ref{lem.det.bordered.det-adj} (applied to }q=i\text{))}}}\nonumber\\
&  \ \ \ \ \ \ \ \ \ \ \left(  \text{here, we have renamed the summation index
}j\text{ as }p\right) \nonumber\\
&  =\sum_{p\in\left\{  1,2,\ldots,n\right\}  }\left(  -1\right)  ^{i+p}%
u_{p}v_{i}\delta_{p,i}\prod_{\substack{j\in\left\{  1,2,\ldots,n\right\}
;\\j\neq p}}d_{j}\nonumber\\
&  =\underbrace{\left(  -1\right)  ^{i+i}}_{\substack{=1\\\text{(since
}i+i=2i\text{ is even)}}}u_{i}v_{i}\underbrace{\delta_{i,i}}%
_{\substack{=1\\\text{(since }i=i\text{)}}}\prod_{\substack{j\in\left\{
1,2,\ldots,n\right\}  ;\\j\neq i}}d_{j}\nonumber\\
&  \ \ \ \ \ \ \ \ \ \ +\sum_{\substack{p\in\left\{  1,2,\ldots,n\right\}
;\\p\neq i}}\left(  -1\right)  ^{i+p}u_{p}v_{i}\underbrace{\delta_{p,i}%
}_{\substack{=0\\\text{(since }p\neq i\text{)}}}\prod_{\substack{j\in\left\{
1,2,\ldots,n\right\}  ;\\j\neq p}}d_{j}\nonumber\\
&  \ \ \ \ \ \ \ \ \ \ \left(
\begin{array}
[c]{c}%
\text{here, we have split off the addend for }p=i\text{ from the sum,}\\
\text{since }i\in\left\{  1,2,\ldots,n\right\}
\end{array}
\right) \nonumber\\
&  =u_{i}v_{i}\prod_{\substack{j\in\left\{  1,2,\ldots,n\right\}  ;\\j\neq
i}}d_{j}+\underbrace{\sum_{\substack{p\in\left\{  1,2,\ldots,n\right\}
;\\p\neq i}}\left(  -1\right)  ^{i+p}u_{p}v_{i}0\prod_{\substack{j\in\left\{
1,2,\ldots,n\right\}  ;\\j\neq p}}d_{j}}_{=0}\nonumber\\
&  =u_{i}v_{i}\prod_{\substack{j\in\left\{  1,2,\ldots,n\right\}  ;\\j\neq
i}}d_{j}. \label{sol.det.bordered.b.3}%
\end{align}
Hence, (\ref{sol.det.bordered.b.2}) becomes%
\begin{align*}
\det\left(
\begin{array}
[c]{cc}%
D & u\\
v & H
\end{array}
\right)   &  =h\underbrace{\det D}_{=d_{1}d_{2}\cdots d_{n}}-\sum_{i=1}%
^{n}\underbrace{\sum_{j=1}^{n}\left(  -1\right)  ^{i+j}u_{j}v_{i}\det\left(
D_{\sim j,\sim i}\right)  }_{\substack{=u_{i}v_{i}\prod_{\substack{j\in
\left\{  1,2,\ldots,n\right\}  ;\\j\neq i}}d_{j}\\\text{(by
(\ref{sol.det.bordered.b.3}))}}}\\
&  =h\cdot\left(  d_{1}d_{2}\cdots d_{n}\right)  -\sum_{i=1}^{n}u_{i}%
v_{i}\prod_{\substack{j\in\left\{  1,2,\ldots,n\right\}  ;\\j\neq i}}d_{j}.
\end{align*}
This solves Exercise \ref{exe.det.bordered} \textbf{(b)}.
\end{proof}

\subsection{\label{sect.sol.det.vdm-pol}Solution to Exercise
\ref{exe.det.vdm-pol}}

Before we start solving Exercise \ref{exe.det.vdm-pol}, we introduce a notation:

\begin{definition}
\label{def.sol.det.vdm-pol.polcoeff}Let $P\in\mathbb{K}\left[  X\right]  $ be
a polynomial. Let $k\in\mathbb{N}$. Then, $\left[  X^{k}\right]  P$ shall
denote the coefficient of $X^{k}$ in $P$.
\end{definition}

For example,%
\[
\left[  X^{2}\right]  \left(  X^{3}+4X^{2}+5X+7\right)
=4\ \ \ \ \ \ \ \ \ \ \text{and}\ \ \ \ \ \ \ \ \ \ \left[  X^{2}\right]
\left(  X+3\right)  =0.
\]

Clearly, each polynomial $P\in\mathbb{K}\left[  X\right]  $ satisfies%
\begin{equation}
P=\sum_{k\in\mathbb{N}}\left(  \left[  X^{k}\right]  P\right)  \cdot X^{k}.
\label{eq.def.sol.det.vdm-pol.polcoeff.1}%
\end{equation}
Moreover, if $n\in\mathbb{N}$, and if $P\in\mathbb{K}\left[  X\right]  $ is a
polynomial satisfying $\deg P\leq n-1$, then%
\begin{equation}
P=\sum_{k=0}^{n-1}\left(  \left[  X^{k}\right]  P\right)  \cdot X^{k}.
\label{eq.def.sol.det.vdm-pol.polcoeff.2}%
\end{equation}

We furthermore ready a simple lemma (the analogue of Exercise \ref{exe.ps4.3}
for upper-triangular matrices):

\begin{lemma}
\label{lem.det.uptriangular}Let $A=\left(  a_{i,j}\right)  _{1\leq i\leq
n,\ 1\leq j\leq n}$ be an $n\times n$-matrix. Assume that $a_{i,j}=0$ for
every $\left(  i,j\right)  \in\left\{  1,2,\ldots,n\right\}  ^{2}$ satisfying
$i>j$. Show that%
\[
\det A=a_{1,1}a_{2,2}\cdots a_{n,n}.
\]

\end{lemma}

\begin{proof}
[Proof of Lemma \ref{lem.det.uptriangular}.]From $A=\left(  a_{i,j}\right)
_{1\leq i\leq n,\ 1\leq j\leq n}$, we obtain $A^{T}=\left(  a_{j,i}\right)
_{1\leq i\leq n,\ 1\leq j\leq n}$ (by the definition of the transpose of a
matrix). Now, recall that%
\begin{equation}
a_{i,j}=0\text{ for every }\left(  i,j\right)  \in\left\{  1,2,\ldots
,n\right\}  ^{2}\text{ satisfying }i>j. \label{pf.lem.det.uptriangular.ass}%
\end{equation}
Hence, $a_{j,i}=0$ for every $\left(  i,j\right)  \in\left\{  1,2,\ldots
,n\right\}  ^{2}$ satisfying $i<j$\ \ \ \ \footnote{\textit{Proof.} Let
$\left(  i,j\right)  \in\left\{  1,2,\ldots,n\right\}  ^{2}$ satisfy $i<j$.
Then, $j>i$ (since $i<j$). Thus, (\ref{pf.lem.det.uptriangular.ass}) (applied
to $\left(  j,i\right)  $ instead of $\left(  i,j\right)  $) yields
$a_{j,i}=0$. Qed.}. Therefore, Exercise \ref{exe.ps4.3} (applied to $A^{T}$
and $a_{j,i}$ instead of $A$ and $a_{i,j}$) yields%
\[
\det\left(  A^{T}\right)  =a_{1,1}a_{2,2}\cdots a_{n,n}%
\]
(since $A^{T}=\left(  a_{j,i}\right)  _{1\leq i\leq n,\ 1\leq j\leq n}$).
However, Exercise \ref{exe.ps4.4} yields $\det\left(  A^{T}\right)  =\det A$.
Comparing these two equalities, we obtain $\det A=a_{1,1}a_{2,2}\cdots
a_{n,n}$. This proves Lemma \ref{lem.det.uptriangular}.
\end{proof}

Now, instead of solving Exercise \ref{exe.det.vdm-pol} directly, we shall
first show the following lemma:

\begin{lemma}
\label{lem.sol.det.vdm-pol.1}Let $n\in\mathbb{N}$. Let $a_{1},a_{2}%
,\ldots,a_{n}\in\mathbb{K}$. Furthermore, for each $j\in\left\{
1,2,\ldots,n\right\}  $, let $P_{j}\in\mathbb{K}\left[  X\right]  $ be a
polynomial such that%
\begin{equation}
\deg\left(  P_{j}\right)  \leq j-1. \label{eq.lem.sol.det.vdm-pol.1.ass}%
\end{equation}
Then,%
\[
\det\left(  \left(  P_{j}\left(  a_{i}\right)  \right)  _{1\leq i\leq
n,\ 1\leq j\leq n}\right)  =\left(  \prod_{j=1}^{n}\left[  X^{j-1}\right]
\left(  P_{j}\right)  \right)  \cdot\det\left(  \left(  a_{i}^{j-1}\right)
_{1\leq i\leq n,\ 1\leq j\leq n}\right)  .
\]

\end{lemma}

\begin{proof}
[Proof of Lemma \ref{lem.sol.det.vdm-pol.1}.]Let $i\in\left\{  1,2,\ldots
,n\right\}  $ and $j\in\left\{  1,2,\ldots,n\right\}  $. The polynomial
$P_{j}$ satisfies $\deg\left(  P_{j}\right)  \leq j-1$ (by
(\ref{eq.lem.sol.det.vdm-pol.1.ass})) and thus $\deg\left(  P_{j}\right)
\leq\underbrace{j}_{\leq n}-1\leq n-1$. Hence,%
\[
P_{j}=\sum_{k=0}^{n-1}\left(  \left[  X^{k}\right]  \left(  P_{j}\right)
\right)  \cdot X^{k}%
\]
(by (\ref{eq.def.sol.det.vdm-pol.polcoeff.2}), applied to $P=P_{j}$).
Substituting $a_{i}$ for $X$ on both sides of this equality, we obtain%
\begin{align}
P_{j}\left(  a_{i}\right)   &  =\sum_{k=0}^{n-1}\underbrace{\left(  \left[
X^{k}\right]  \left(  P_{j}\right)  \right)  \cdot a_{i}^{k}}_{=a_{i}^{k}%
\cdot\left[  X^{k}\right]  \left(  P_{j}\right)  }=\sum_{k=0}^{n-1}a_{i}%
^{k}\cdot\left[  X^{k}\right]  \left(  P_{j}\right) \nonumber\\
&  =\sum_{k=1}^{n}a_{i}^{k-1}\cdot\left[  X^{k-1}\right]  \left(
P_{j}\right)  \label{pf.lem.sol.det.vdm-pol.1.1}%
\end{align}
(here, we have substituted $k-1$ for $k$ in the sum).

Forget that we fixed $i$ and $j$. We thus have proved the equality
(\ref{pf.lem.sol.det.vdm-pol.1.1}) for all $i\in\left\{  1,2,\ldots,n\right\}
$ and $j\in\left\{  1,2,\ldots,n\right\}  $. In other words, we have%
\begin{equation}
\left(  P_{j}\left(  a_{i}\right)  \right)  _{1\leq i\leq n,\ 1\leq j\leq
n}=\left(  \sum_{k=1}^{n}a_{i}^{k-1}\cdot\left[  X^{k-1}\right]  \left(
P_{j}\right)  \right)  _{1\leq i\leq n,\ 1\leq j\leq n}.
\label{pf.lem.sol.det.vdm-pol.1.2}%
\end{equation}

Define an $n\times n$-matrix $B\in\mathbb{K}^{n\times n}$ by
\begin{equation}
B=\left(  \left[  X^{i-1}\right]  \left(  P_{j}\right)  \right)  _{1\leq i\leq
n,\ 1\leq j\leq n}. \label{pf.lem.sol.det.vdm-pol.1.B}%
\end{equation}

We have $\left[  X^{i-1}\right]  \left(  P_{j}\right)  =0$ for every $\left(
i,j\right)  \in\left\{  1,2,\ldots,n\right\}  ^{2}$ satisfying $i>j$%
\ \ \ \ \footnote{\textit{Proof.} Let $\left(  i,j\right)  \in\left\{
1,2,\ldots,n\right\}  ^{2}$ be such that $i>j$. From $i>j$, we obtain $i\geq
j+1$ (since $i$ and $j$ are integers), so that $i-1\geq j$.
\par
We have $\deg\left(  P_{j}\right)  \leq j-1$ (by
(\ref{eq.lem.sol.det.vdm-pol.1.ass})). In other words, the polynomial $P_{j}$
has degree $\leq j-1$. Hence, the coefficients of the monomials $X^{j}%
,X^{j+1},X^{j+2},\ldots$ in this polynomial $P_{j}$ are $0$. In other words,
if $k$ is any integer satisfying $k\geq j$, then the coefficient of the
monomial $X^{k}$ in the polynomial $P_{j}$ is $0$. Applying this to $k=i-1$,
we obtain that the coefficient of the monomial $X^{i-1}$ in the polynomial
$P_{j}$ is $0$ (since $i-1\geq j$). In other words, $\left[  X^{i-1}\right]
\left(  P_{j}\right)  =0$ (since $\left[  X^{i-1}\right]  \left(
P_{j}\right)  $ was defined to be the coefficient of the monomial $X^{i-1}$ in
the polynomial $P_{j}$). Qed.}. Hence, Lemma \ref{lem.det.uptriangular}
(applied to $B$ and $\left[  X^{i-1}\right]  \left(  P_{j}\right)  $ instead
of $A$ and $a_{i,j}$) yields
\begin{align}
\det B  &  =\left(  \left[  X^{1-1}\right]  \left(  P_{1}\right)  \right)
\left(  \left[  X^{2-1}\right]  \left(  P_{2}\right)  \right)  \cdots\left(
\left[  X^{n-1}\right]  \left(  P_{n}\right)  \right) \nonumber\\
&  =\prod_{j=1}^{n}\left[  X^{j-1}\right]  \left(  P_{j}\right)  .
\label{pf.lem.sol.det.vdm-pol.1.detB=}%
\end{align}

Define an $n\times n$-matrix $A\in\mathbb{K}^{n\times n}$ by%
\begin{equation}
A=\left(  a_{i}^{j-1}\right)  _{1\leq i\leq n,\ 1\leq j\leq n}.
\label{pf.lem.sol.det.vdm-pol.1.A}%
\end{equation}

The definition of the product of two matrices yields%
\[
AB=\left(  \sum_{k=1}^{n}a_{i}^{k-1}\cdot\left[  X^{k-1}\right]  \left(
P_{j}\right)  \right)  _{1\leq i\leq n,\ 1\leq j\leq n}%
\]
(because of (\ref{pf.lem.sol.det.vdm-pol.1.A}) and
(\ref{pf.lem.sol.det.vdm-pol.1.B})). Comparing this equality with
(\ref{pf.lem.sol.det.vdm-pol.1.2}), we obtain%
\[
\left(  P_{j}\left(  a_{i}\right)  \right)  _{1\leq i\leq n,\ 1\leq j\leq
n}=AB.
\]
Hence,%
\begin{align*}
\det\left(  \underbrace{\left(  P_{j}\left(  a_{i}\right)  \right)  _{1\leq
i\leq n,\ 1\leq j\leq n}}_{=AB}\right)   &  =\det\left(  AB\right)  =\det
A\cdot\det B\ \ \ \ \ \ \ \ \ \ \left(  \text{by Theorem \ref{thm.det(AB)}%
}\right) \\
&  =\underbrace{\det B}_{\substack{=\prod_{j=1}^{n}\left[  X^{j-1}\right]
\left(  P_{j}\right)  \\\text{(by (\ref{pf.lem.sol.det.vdm-pol.1.detB=}))}%
}}\cdot\det\underbrace{A}_{=\left(  a_{i}^{j-1}\right)  _{1\leq i\leq
n,\ 1\leq j\leq n}}\\
&  =\left(  \prod_{j=1}^{n}\left[  X^{j-1}\right]  \left(  P_{j}\right)
\right)  \cdot\det\left(  \left(  a_{i}^{j-1}\right)  _{1\leq i\leq n,\ 1\leq
j\leq n}\right)  .
\end{align*}
This proves Lemma \ref{lem.sol.det.vdm-pol.1}.
\end{proof}

We could now easily solve Exercise \ref{exe.det.vdm-pol} by combining Lemma
\ref{lem.sol.det.vdm-pol.1} with Theorem \ref{thm.vander-det} \textbf{(c)}.
However, let us extract some additional usefulness from Lemma
\ref{lem.sol.det.vdm-pol.1} by proving Theorem \ref{thm.vander-det}
\textbf{(c)} again:

\begin{proof}
[Third proof of Theorem \ref{thm.vander-det} \textbf{(c)}.]Let us rename the
$n$ elements $x_{1},x_{2},\ldots,x_{n}$ of $\mathbb{K}$ as $a_{1},a_{2}%
,\ldots,a_{n}$. Thus, we must prove that%
\begin{equation}
\det\left(  \left(  a_{i}^{j-1}\right)  _{1\leq i\leq n,\ 1\leq j\leq
n}\right)  =\prod_{1\leq j<i\leq n}\left(  a_{i}-a_{j}\right)  .
\label{pf.thm.vander-det.c.3rd.goal}%
\end{equation}

For each $j\in\left\{  1,2,\ldots,n\right\}  $, we define a polynomial
$P_{j}\in\mathbb{K}\left[  X\right]  $ by%
\[
P_{j}=\left(  X-a_{1}\right)  \left(  X-a_{2}\right)  \cdots\left(
X-a_{j-1}\right)  .
\]
This polynomial $P_{j}$ is a product of $j-1$ monic polynomials of degree $1$
(namely, of the $j-1$ polynomials $X-a_{1},X-a_{2},\ldots,X-a_{j-1}$), and
thus itself is a monic polynomial of degree $\underbrace{1+1+\cdots
+1}_{j-1\text{ times}}$ (since a product of several monic polynomials is
always a monic polynomial, and its degree is the sum of their degrees). In
other words, $P_{j}$ is a monic polynomial of degree $j-1$ (since
$\underbrace{1+1+\cdots+1}_{j-1\text{ times}}=\left(  j-1\right)  \cdot
1=j-1$). Hence, in particular, this polynomial $P_{j}$ has degree $j-1$, so
that $\deg\left(  P_{j}\right)  =j-1\leq j-1$. Thus, Lemma
\ref{lem.sol.det.vdm-pol.1} yields
\[
\det\left(  \left(  P_{j}\left(  a_{i}\right)  \right)  _{1\leq i\leq
n,\ 1\leq j\leq n}\right)  =\left(  \prod_{j=1}^{n}\left[  X^{j-1}\right]
\left(  P_{j}\right)  \right)  \cdot\det\left(  \left(  a_{i}^{j-1}\right)
_{1\leq i\leq n,\ 1\leq j\leq n}\right)  .
\]

\begin{vershort}
However, for each $j\in\left\{  1,2,\ldots,n\right\}  $, we have $\left[
X^{j-1}\right]  \left(  P_{j}\right)  =1$ (since we have just shown that
$P_{j}$ is a monic polynomial of degree $j-1$). Multiplying these equalities
over all $j\in\left\{  1,2,\ldots,n\right\}  $, we obtain $\prod_{j=1}%
^{n}\left[  X^{j-1}\right]  \left(  P_{j}\right)  =\prod_{j=1}^{n}1=1$.
\end{vershort}

\begin{verlong}
However, for each $j\in\left\{  1,2,\ldots,n\right\}  $, we have $\left[
X^{j-1}\right]  \left(  P_{j}\right)  =1$\ \ \ \ \footnote{\textit{Proof.} Let
$j\in\left\{  1,2,\ldots,n\right\}  $. Then, $P_{j}$ is a monic polynomial of
degree $j-1$ (as we have seen above). Hence, the coefficient of $X^{j-1}$ in
this polynomial $P_{j}$ is $1$. In other words, $\left[  X^{j-1}\right]
\left(  P_{j}\right)  =1$ (since $\left[  X^{j-1}\right]  \left(
P_{j}\right)  $ is defined to be the coefficient of $X^{j-1}$ in this
polynomial $P_{j}$). Qed.}. Multiplying these equalities over all
$j\in\left\{  1,2,\ldots,n\right\}  $, we obtain $\prod_{j=1}^{n}\left[
X^{j-1}\right]  \left(  P_{j}\right)  =\prod_{j=1}^{n}1=1$.
\end{verlong}

Thus,%
\begin{align}
\det\left(  \left(  P_{j}\left(  a_{i}\right)  \right)  _{1\leq i\leq
n,\ 1\leq j\leq n}\right)   &  =\underbrace{\left(  \prod_{j=1}^{n}\left[
X^{j-1}\right]  \left(  P_{j}\right)  \right)  }_{=1}\cdot\det\left(  \left(
a_{i}^{j-1}\right)  _{1\leq i\leq n,\ 1\leq j\leq n}\right) \nonumber\\
&  =\det\left(  \left(  a_{i}^{j-1}\right)  _{1\leq i\leq n,\ 1\leq j\leq
n}\right)  . \label{pf.thm.vander-det.c.3rd.4}%
\end{align}

Note that $P_{j}\left(  a_{i}\right)  =0$ for every $\left(  i,j\right)
\in\left\{  1,2,\ldots,n\right\}  ^{2}$ satisfying $i<j$%
\ \ \ \ \footnote{\textit{Proof.} Let $\left(  i,j\right)  \in\left\{
1,2,\ldots,n\right\}  ^{2}$ be such that $i<j$. From $i<j$, we obtain $i\leq
j-1$ (since $i$ and $j$ are integers), so that $i\in\left\{  1,2,\ldots
,j-1\right\}  $.
\par
However, the definition of $P_{j}$ yields $P_{j}=\left(  X-a_{1}\right)
\left(  X-a_{2}\right)  \cdots\left(  X-a_{j-1}\right)  =\prod_{k\in\left\{
1,2,\ldots,j-1\right\}  }\left(  X-a_{k}\right)  $. Substituting $a_{i}$ for
$X$ on both sides of this equality, we obtain%
\begin{align*}
P_{j}\left(  a_{i}\right)   &  =\prod_{k\in\left\{  1,2,\ldots,j-1\right\}
}\left(  a_{i}-a_{k}\right)  =\underbrace{\left(  a_{i}-a_{i}\right)  }%
_{=0}\cdot\prod_{\substack{k\in\left\{  1,2,\ldots,j-1\right\}  ;\\k\neq
i}}\left(  a_{i}-a_{k}\right) \\
&  \ \ \ \ \ \ \ \ \ \ \ \ \ \ \ \ \ \ \ \ \left(
\begin{array}
[c]{c}%
\text{here, we have split off the factor for }k=i\\
\text{from the product (since }i\in\left\{  1,2,\ldots,j-1\right\}  \text{)}%
\end{array}
\right) \\
&  =0,
\end{align*}
qed.}. Hence, Exercise \ref{exe.ps4.3} (applied to $\left(  P_{j}\left(
a_{i}\right)  \right)  _{1\leq i\leq n,\ 1\leq j\leq n}$ and $P_{j}\left(
a_{i}\right)  $ instead of $A$ and $a_{i,j}$) yields%
\[
\det\left(  \left(  P_{j}\left(  a_{i}\right)  \right)  _{1\leq i\leq
n,\ 1\leq j\leq n}\right)  =P_{1}\left(  a_{1}\right)  \cdot P_{2}\left(
a_{2}\right)  \cdot\cdots\cdot P_{n}\left(  a_{n}\right)  =\prod_{j=1}%
^{n}P_{j}\left(  a_{j}\right)  .
\]

However, each $j\in\left\{  1,2,\ldots,n\right\}  $ satisfies $P_{j}\left(
a_{j}\right)  =\prod_{i=1}^{j-1}\left(  a_{j}-a_{i}\right)  $%
\ \ \ \ \footnote{\textit{Proof.} Let $j\in\left\{  1,2,\ldots,n\right\}  $.
The definition of $P_{j}$ yields $P_{j}=\left(  X-a_{1}\right)  \left(
X-a_{2}\right)  \cdots\left(  X-a_{j-1}\right)  =\prod_{i=1}^{j-1}\left(
X-a_{i}\right)  $. Substituting $a_{j}$ for $X$ on both sides of this
equality, we obtain $P_{j}\left(  a_{j}\right)  =\prod_{i=1}^{j-1}\left(
a_{j}-a_{i}\right)  $, qed.}. Multiplying these equalities over all
$j\in\left\{  1,2,\ldots,n\right\}  $, we obtain%
\[
\prod_{j=1}^{n}P_{j}\left(  a_{j}\right)  =\underbrace{\prod_{j=1}%
^{n}\ \ \prod_{i=1}^{j-1}}_{=\prod_{1\leq i<j\leq n}}\left(  a_{j}%
-a_{i}\right)  =\prod_{1\leq i<j\leq n}\left(  a_{j}-a_{i}\right)
=\prod_{1\leq j<i\leq n}\left(  a_{i}-a_{j}\right)
\]
(here, we have renamed the index $\left(  i,j\right)  $ as $\left(
j,i\right)  $ in the product). Now, from (\ref{pf.thm.vander-det.c.3rd.4}), we
obtain%
\begin{align*}
\det\left(  \left(  a_{i}^{j-1}\right)  _{1\leq i\leq n,\ 1\leq j\leq
n}\right)   &  =\det\left(  \left(  P_{j}\left(  a_{i}\right)  \right)
_{1\leq i\leq n,\ 1\leq j\leq n}\right)  =\prod_{j=1}^{n}P_{j}\left(
a_{j}\right) \\
&  =\prod_{1\leq j<i\leq n}\left(  a_{i}-a_{j}\right)  .
\end{align*}
Hence, we have proved (\ref{pf.thm.vander-det.c.3rd.goal}). In other words, we
have proved Theorem \ref{thm.vander-det} \textbf{(c)} once again (since our
$a_{1},a_{2},\ldots,a_{n}$ are just the renamed $x_{1},x_{2},\ldots,x_{n}$).
\end{proof}

Let us now solve Exercise \ref{exe.det.vdm-pol}:

\begin{vershort}
\begin{proof}
[Solution to Exercise \ref{exe.det.vdm-pol}.]We have%
\begin{equation}
\left[  X^{j-1}\right]  \left(  P_{j}\right)  =c_{j}%
\ \ \ \ \ \ \ \ \ \ \text{for each }j\in\left\{  1,2,\ldots,n\right\}
\label{sol.det.vdm-pol.short.=cj}%
\end{equation}
(since both $\left[  X^{j-1}\right]  \left(  P_{j}\right)  $ and $c_{j}$ are
defined to be the coefficient of $X^{j-1}$ in the polynomial $P_{j}$). Now,
Lemma \ref{lem.sol.det.vdm-pol.1} yields%
\begin{align*}
\det\left(  \left(  P_{j}\left(  a_{i}\right)  \right)  _{1\leq i\leq
n,\ 1\leq j\leq n}\right)   &  =\left(  \prod_{j=1}^{n}\underbrace{\left[
X^{j-1}\right]  \left(  P_{j}\right)  }_{\substack{=c_{j}\\\text{(by
(\ref{sol.det.vdm-pol.short.=cj}))}}}\right)  \cdot\underbrace{\det\left(
\left(  a_{i}^{j-1}\right)  _{1\leq i\leq n,\ 1\leq j\leq n}\right)
}_{\substack{=\prod_{1\leq j<i\leq n}\left(  a_{i}-a_{j}\right)  \\\text{(by
Theorem \ref{thm.vander-det} \textbf{(c)},}\\\text{applied to }x_{i}%
=a_{i}\text{)}}}\\
&  =\left(  \prod_{j=1}^{n}c_{j}\right)  \cdot\prod_{1\leq j<i\leq n}\left(
a_{i}-a_{j}\right)  .
\end{align*}
This solves Exercise \ref{exe.det.vdm-pol}.
\end{proof}
\end{vershort}

\begin{verlong}
\begin{proof}
[Solution to Exercise \ref{exe.det.vdm-pol}.]We have%
\begin{equation}
\left[  X^{j-1}\right]  \left(  P_{j}\right)  =c_{j}%
\ \ \ \ \ \ \ \ \ \ \text{for each }j\in\left\{  1,2,\ldots,n\right\}  .
\label{sol.det.vdm-pol.=cj}%
\end{equation}
\footnote{\textit{Proof of (\ref{sol.det.vdm-pol.=cj}):} Let $j\in\left\{
1,2,\ldots,n\right\}  $. Then,%
\[
c_{j}=\left(  \text{the coefficient of }X^{j-1}\text{ in the polynomial }%
P_{j}\right)
\]
(by the definition of $c_{j}$). On the other hand,
\[
\left[  X^{j-1}\right]  \left(  P_{j}\right)  =\left(  \text{the coefficient
of }X^{j-1}\text{ in the polynomial }P_{j}\right)
\]
(by the definition of $\left[  X^{j-1}\right]  \left(  P_{j}\right)  $).
Comparing these two equalities, we obtain $\left[  X^{j-1}\right]  \left(
P_{j}\right)  =c_{j}$. This proves (\ref{sol.det.vdm-pol.=cj}).}

Now, Lemma \ref{lem.sol.det.vdm-pol.1} yields%
\begin{align*}
\det\left(  \left(  P_{j}\left(  a_{i}\right)  \right)  _{1\leq i\leq
n,\ 1\leq j\leq n}\right)   &  =\left(  \prod_{j=1}^{n}\underbrace{\left[
X^{j-1}\right]  \left(  P_{j}\right)  }_{\substack{=c_{j}\\\text{(by
(\ref{sol.det.vdm-pol.=cj}))}}}\right)  \cdot\underbrace{\det\left(  \left(
a_{i}^{j-1}\right)  _{1\leq i\leq n,\ 1\leq j\leq n}\right)  }%
_{\substack{=\prod_{1\leq j<i\leq n}\left(  a_{i}-a_{j}\right)  \\\text{(by
Theorem \ref{thm.vander-det} \textbf{(c)},}\\\text{applied to }x_{i}%
=a_{i}\text{)}}}\\
&  =\left(  \prod_{j=1}^{n}c_{j}\right)  \cdot\prod_{1\leq j<i\leq n}\left(
a_{i}-a_{j}\right)  .
\end{align*}
This solves Exercise \ref{exe.det.vdm-pol}.
\end{proof}
\end{verlong}

\subsection{Solution to Exercise \ref{exe.det.kratt-lem6}}

We will be using Definition \ref{def.sol.det.vdm-pol.polcoeff} throughout this
section. The following two facts about polynomials are obvious:

\begin{lemma}
\label{lem.sol.det.kratt-lem6.pol-lem1}Let $n\in\mathbb{N}$. Let $a$ and $b$
be two polynomials in $\mathbb{K}\left[  X\right]  $ that satisfy $\deg a\leq
n$ and $\deg b\leq n$. Then, $\deg\left(  a-b\right)  \leq n$.
\end{lemma}

\begin{verlong}
\begin{proof}
[Proof of Lemma \ref{lem.sol.det.kratt-lem6.pol-div.1}.]If $m$ is any integer
satisfying $m>n$, then
\[
\left[  X^{m}\right]  \left(  a-b\right)  =\underbrace{\left[  X^{m}\right]
a}_{\substack{=0\\\text{(since }m>n\geq\deg a\text{)}}}-\underbrace{\left[
X^{m}\right]  b}_{\substack{=0\\\text{(since }m>n\geq\deg b\text{)}}}=0.
\]
In other words, if $m$ is any integer satisfying $m>n$, then the coefficient
of the monomial $X^{m}$ in $a-b$ is $0$. In other words, all the monomials
$X^{n+1},X^{n+2},X^{n+3},\ldots$ appear with coefficient $0$ in the polynomial
$a-b$. Thus, this polynomial $a-b$ has degree $\leq n$. In other words,
$\deg\left(  a-b\right)  \leq n$. This proves Lemma
\ref{lem.sol.det.kratt-lem6.pol-div.1}.
\end{proof}
\end{verlong}

\begin{lemma}
\label{lem.sol.det.kratt-lem6.pol-lem2}Let $n\in\mathbb{N}$. Let $a$ be a
polynomial in $\mathbb{K}\left[  X\right]  $ that satisfies $\deg a\leq n$ and
$\left[  X^{n}\right]  a=0$. Then, $\deg a\leq n-1$.
\end{lemma}

\begin{proof}
[Proof of Lemma \ref{lem.sol.det.kratt-lem6.pol-lem2}.]If we had $\deg a=n$,
then we would have $\left[  X^{n}\right]  a\neq0$, which would contradict
$\left[  X^{n}\right]  a=0$. Thus, we cannot have $\deg a=n$. Therefore, $\deg
a\neq n$. Combining this with $\deg a\leq n$, we obtain $\deg a<n$, so that
$\deg a\leq n-1$ (since $\deg a$ and $n$ are either integers or $-\infty$).
This proves Lemma \ref{lem.sol.det.kratt-lem6.pol-lem2}.
\end{proof}

Before solving Exercise \ref{exe.det.kratt-lem6}, we will state a lemma which
generalizes the obvious fact that a polynomial of degree $\leq j-1$ in
$\mathbb{K}\left[  X\right]  $ can be written as $d_{1}X^{j-1}+d_{2}%
X^{j-2}+\cdots+d_{j}X^{j-j}$ for some $d_{1},d_{2},\ldots,d_{j}\in\mathbb{K}$:

\begin{lemma}
\label{lem.sol.det.kratt-lem6.pol-div.1}Let $j\in\mathbb{N}$. For each
$k\in\left\{  1,2,\ldots,j\right\}  $, let $R_{k}$ be a monic polynomial of
degree $j-k$ in $\mathbb{K}\left[  X\right]  $.

Let $P\in\mathbb{K}\left[  X\right]  $ be a polynomial such that $\deg P\leq
j-1$. Then, there exist $j$ elements $d_{1},d_{2},\ldots,d_{j}$ of
$\mathbb{K}$ such that%
\[
P=d_{1}R_{1}+d_{2}R_{2}+\cdots+d_{j}R_{j}.
\]

\end{lemma}

\begin{proof}
[Proof of Lemma \ref{lem.sol.det.kratt-lem6.pol-div.1}.]For each $k\in\left\{
1,2,\ldots,j\right\}  $, we have%
\begin{equation}
\left[  X^{j-k}\right]  \left(  R_{k}\right)  =1.
\label{pf.lem.sol.det.kratt-lem6.pol-div.1.1}%
\end{equation}
\footnote{\textit{Proof of (\ref{pf.lem.sol.det.kratt-lem6.pol-div.1.1}):} Let
$k\in\left\{  1,2,\ldots,j\right\}  $. Then, $R_{k}$ is a monic polynomial of
degree $j-k$ (by the definition of $R_{k}$). Hence, the coefficient of
$X^{j-k}$ in this polynomial $R_{k}$ is $1$. In other words, we have $\left[
X^{j-k}\right]  \left(  R_{k}\right)  =1$ (since $\left[  X^{j-k}\right]
\left(  R_{k}\right)  $ is defined to be the coefficient of $X^{j-k}$ in this
polynomial $R_{k}$). This proves (\ref{pf.lem.sol.det.kratt-lem6.pol-div.1.1}%
).}

We now define $j$ elements $c_{1},c_{2},\ldots,c_{j}$ of $\mathbb{K}$
recursively as follows: For each $m\in\left\{  1,2,\ldots,j\right\}  $, we
assume that the first $m-1$ elements $c_{1},c_{2},\ldots,c_{m-1}$ have already
been defined, and we define $c_{m}$ by%
\begin{equation}
c_{m}=\left[  X^{j-m}\right]  \left(  P-\sum_{k=1}^{m-1}c_{k}R_{k}\right)  .
\label{pf.lem.sol.det.kratt-lem6.pol-div.1.di=}%
\end{equation}
Thus, we have defined $j$ elements $c_{1},c_{2},\ldots,c_{j}$ of $\mathbb{K}$.

Now, we shall show that
\begin{equation}
\deg\left(  P-\sum_{k=1}^{i}c_{k}R_{k}\right)  \leq
j-i-1\ \ \ \ \ \ \ \ \ \ \text{for each }i\in\left\{  0,1,\ldots,j\right\}  .
\label{pf.lem.sol.det.kratt-lem6.pol-div.1.deg}%
\end{equation}

[\textit{Proof of (\ref{pf.lem.sol.det.kratt-lem6.pol-div.1.deg}):} We shall
prove (\ref{pf.lem.sol.det.kratt-lem6.pol-div.1.deg}) by induction on $i$:

\textit{Induction base:} It is easy to see that
(\ref{pf.lem.sol.det.kratt-lem6.pol-div.1.deg}) holds for $i=0$%
\ \ \ \ \footnote{\textit{Proof.} We have $P-\underbrace{\sum_{k=1}^{0}%
c_{k}R_{k}}_{=\left(  \text{empty sum}\right)  =0}=P$ and thus $\deg\left(
P-\sum_{k=1}^{0}c_{k}R_{k}\right)  =\deg P\leq j-1=j-0-1$. In other words,
(\ref{pf.lem.sol.det.kratt-lem6.pol-div.1.deg}) holds for $i=0$.}. This
completes the induction base.

\textit{Induction step:} Let $m\in\left\{  1,2,\ldots,j\right\}  $. Assume (as
the induction hypothesis) that (\ref{pf.lem.sol.det.kratt-lem6.pol-div.1.deg})
holds for $i=m-1$. We shall show that
(\ref{pf.lem.sol.det.kratt-lem6.pol-div.1.deg}) holds for $i=m$.

We have assumed that (\ref{pf.lem.sol.det.kratt-lem6.pol-div.1.deg}) holds for
$i=m-1$. In other words, we have%
\[
\deg\left(  P-\sum_{k=1}^{m-1}c_{k}R_{k}\right)  \leq j-\left(  m-1\right)
-1.
\]

\begin{vershort}
However, $R_{m}$ is a monic polynomial of degree $j-m$ (by the definition of
$R_{m}$). Thus, in particular, $\deg\left(  R_{m}\right)  \leq j-m$. Hence,
$\deg\left(  c_{m}R_{m}\right)  \leq j-m$ as well (because multiplying a
polynomial by a constant will never increase the degree of the polynomial).
\end{vershort}

\begin{verlong}
Now, recall that $R_{k}$ is a monic polynomial of degree $j-k$ whenever
$k\in\left\{  1,2,\ldots,j\right\}  $. Applying this to $k=m$, we conclude
that $R_{m}$ is a monic polynomial of degree $j-m$. Hence, $R_{m}$ has degree
$j-m$. In other words, $\deg\left(  R_{m}\right)  =j-m$, so that $\deg\left(
R_{m}\right)  \leq j-m$. However, if $a\in\mathbb{K}\left[  X\right]  $ is any
polynomial and $c\in\mathbb{K}$ is any constant, then $\deg\left(  ca\right)
\leq\deg a$. Applying this to $c=c_{m}$ and $a=R_{m}$, we obtain $\deg\left(
c_{m}R_{m}\right)  \leq\deg\left(  R_{m}\right)  \leq j-m$.

Furthermore, $m\in\left\{  1,2,\ldots,j\right\}  $ and thus $m\leq j$. Hence,
$j-m\geq0$, so that $j-m\in\mathbb{N}$.
\end{verlong}

Now, Lemma \ref{lem.sol.det.kratt-lem6.pol-lem1} (applied to $a=P-\sum
_{k=1}^{m-1}c_{k}R_{k}$ and $b=c_{m}R_{m}$ and $n=j-m$) yields%
\[
\deg\left(  \left(  P-\sum_{k=1}^{m-1}c_{k}R_{k}\right)  -c_{m}R_{m}\right)
\leq j-m
\]
(since $\deg\left(  P-\sum_{k=1}^{m-1}c_{k}R_{k}\right)  \leq j-\left(
m-1\right)  -1=j-m$ and $\deg\left(  c_{m}R_{m}\right)  \leq j-m$). In view of%
\begin{align}
P-\underbrace{\sum_{k=1}^{m}c_{k}R_{k}}_{=\sum_{k=1}^{m-1}c_{k}R_{k}%
+c_{m}R_{m}}  &  =P-\left(  \sum_{k=1}^{m-1}c_{k}R_{k}+c_{m}R_{m}\right)
\nonumber\\
&  =\left(  P-\sum_{k=1}^{m-1}c_{k}R_{k}\right)  -c_{m}R_{m},
\label{pf.lem.sol.det.kratt-lem6.pol-div.1.deg.pf.1}%
\end{align}
we can rewrite this as%
\[
\deg\left(  P-\sum_{k=1}^{m}c_{k}R_{k}\right)  \leq j-m.
\]

Now,%
\begin{align*}
&  \left[  X^{j-m}\right]  \underbrace{\left(  P-\sum_{k=1}^{m}c_{k}%
R_{k}\right)  }_{\substack{=\left(  P-\sum_{k=1}^{m-1}c_{k}R_{k}\right)
-c_{m}R_{m}\\\text{(by (\ref{pf.lem.sol.det.kratt-lem6.pol-div.1.deg.pf.1}))}%
}}\\
&  =\left[  X^{j-m}\right]  \left(  \left(  P-\sum_{k=1}^{m-1}c_{k}%
R_{k}\right)  -c_{m}R_{m}\right) \\
&  =\left[  X^{j-m}\right]  \left(  P-\sum_{k=1}^{m-1}c_{k}R_{k}\right)
-\underbrace{\left[  X^{j-m}\right]  \left(  c_{m}R_{m}\right)  }%
_{\substack{=c_{m}\cdot\left[  X^{j-m}\right]  \left(  R_{m}\right)
\\\text{(since }\left[  X^{j-m}\right]  \left(  \lambda a\right)
=\lambda\cdot\left[  X^{j-m}\right]  a\\\text{for any }\lambda\in
\mathbb{K}\text{ and any }a\in\mathbb{K}\left[  X\right]  \text{)}}}\\
&  \ \ \ \ \ \ \ \ \ \ \ \ \ \ \ \ \ \ \ \ \left(
\begin{array}
[c]{c}%
\text{since }\left[  X^{j-m}\right]  \left(  a-b\right)  =\left[
X^{j-m}\right]  a-\left[  X^{j-m}\right]  b\\
\text{for any two polynomials }a\text{ and }b\text{ in }\mathbb{K}\left[
X\right]
\end{array}
\right) \\
&  =\underbrace{\left[  X^{j-m}\right]  \left(  P-\sum_{k=1}^{m-1}c_{k}%
R_{k}\right)  }_{\substack{=c_{m}\\\text{(by
(\ref{pf.lem.sol.det.kratt-lem6.pol-div.1.di=}))}}}-c_{m}\cdot
\underbrace{\left[  X^{j-m}\right]  \left(  R_{m}\right)  }%
_{\substack{=1\\\text{(by (\ref{pf.lem.sol.det.kratt-lem6.pol-div.1.1}%
),}\\\text{applied to }k=m\text{)}}}\\
&  =c_{m}-c_{m}=0.
\end{align*}
Hence, Lemma \ref{lem.sol.det.kratt-lem6.pol-lem2} (applied to $n=j-m$ and
$a=P-\sum_{k=1}^{m}c_{k}R_{k}$) yields
\[
\deg\left(  P-\sum_{k=1}^{m}c_{k}R_{k}\right)  \leq j-m-1
\]
(since we also know that $\deg\left(  P-\sum_{k=1}^{m}c_{k}R_{k}\right)  \leq
j-m$). In other words, (\ref{pf.lem.sol.det.kratt-lem6.pol-div.1.deg}) holds
for $i=m$. This completes the induction step. Thus,
(\ref{pf.lem.sol.det.kratt-lem6.pol-div.1.deg}) is proved by induction.]

Now, (\ref{pf.lem.sol.det.kratt-lem6.pol-div.1.deg}) (applied to $i=j$) yields%
\[
\deg\left(  P-\sum_{k=1}^{j}c_{k}R_{k}\right)  \leq j-j-1=-1<0.
\]
In other words, the polynomial $P-\sum_{k=1}^{j}c_{k}R_{k}$ has negative
degree. Therefore, this polynomial $P-\sum_{k=1}^{j}c_{k}R_{k}$ must be $0$
(since the only polynomial that has negative degree is $0$). In other words,
$P-\sum_{k=1}^{j}c_{k}R_{k}=0$, so that $P=\sum_{k=1}^{j}c_{k}R_{k}=c_{1}%
R_{1}+c_{2}R_{2}+\cdots+c_{j}R_{j}$.

Hence, there exist $j$ elements $d_{1},d_{2},\ldots,d_{j}$ of $\mathbb{K}$
such that $P=d_{1}R_{1}+d_{2}R_{2}+\cdots+d_{j}R_{j}$ (namely, $d_{i}=c_{i}$).
This proves Lemma \ref{lem.sol.det.kratt-lem6.pol-div.1}.
\end{proof}

Our next lemma is an analogue for Lemma \ref{lem.vander-det.lem-rearr} for
columns instead of rows:

\begin{lemma}
\label{lem.vander-det.lem-rearr-col}Let $n\in\mathbb{N}$. Let $\left(
a_{i,j}\right)  _{1\leq i\leq n,\ 1\leq j\leq n}$ be an $n\times n$-matrix.
Then,%
\[
\det\left(  \left(  a_{i,n+1-j}\right)  _{1\leq i\leq n,\ 1\leq j\leq
n}\right)  =\left(  -1\right)  ^{n\left(  n-1\right)  /2}\det\left(  \left(
a_{i,j}\right)  _{1\leq i\leq n,\ 1\leq j\leq n}\right)  .
\]

\end{lemma}

\begin{proof}
[Proof of Lemma \ref{lem.vander-det.lem-rearr-col}.]Lemma
\ref{lem.vander-det.lem-rearr} (applied to $a_{j,i}$ instead of $a_{i,j}$)
yields%
\begin{align}
&  \det\left(  \left(  a_{j,n+1-i}\right)  _{1\leq i\leq n,\ 1\leq j\leq
n}\right) \nonumber\\
&  =\left(  -1\right)  ^{n\left(  n-1\right)  /2}\det\left(  \left(
a_{j,i}\right)  _{1\leq i\leq n,\ 1\leq j\leq n}\right)  .
\label{pf.lem.vander-det.lem-rearr-col.1}%
\end{align}

However, we have $\left(  \left(  a_{i,j}\right)  _{1\leq i\leq n,\ 1\leq
j\leq n}\right)  ^{T}=\left(  a_{j,i}\right)  _{1\leq i\leq n,\ 1\leq j\leq
n}$ (by the definition of the transpose of a matrix) and therefore%
\[
\det\underbrace{\left(  \left(  \left(  a_{i,j}\right)  _{1\leq i\leq
n,\ 1\leq j\leq n}\right)  ^{T}\right)  }_{=\left(  a_{j,i}\right)  _{1\leq
i\leq n,\ 1\leq j\leq n}}=\det\left(  \left(  a_{j,i}\right)  _{1\leq i\leq
n,\ 1\leq j\leq n}\right)  ,
\]
so that%
\begin{align}
\det\left(  \left(  a_{j,i}\right)  _{1\leq i\leq n,\ 1\leq j\leq n}\right)
&  =\det\left(  \left(  \left(  a_{i,j}\right)  _{1\leq i\leq n,\ 1\leq j\leq
n}\right)  ^{T}\right) \nonumber\\
&  =\det\left(  \left(  a_{i,j}\right)  _{1\leq i\leq n,\ 1\leq j\leq
n}\right)  \label{pf.lem.vander-det.lem-rearr-col.2}%
\end{align}
(by Exercise \ref{exe.ps4.4}, applied to $A=\left(  a_{i,j}\right)  _{1\leq
i\leq n,\ 1\leq j\leq n}$).

Furthermore, we have $\left(  \left(  a_{i,n+1-j}\right)  _{1\leq i\leq
n,\ 1\leq j\leq n}\right)  ^{T}=\left(  a_{j,n+1-i}\right)  _{1\leq i\leq
n,\ 1\leq j\leq n}$ (by the definition of the transpose of a matrix) and
therefore%
\[
\det\underbrace{\left(  \left(  \left(  a_{i,n+1-j}\right)  _{1\leq i\leq
n,\ 1\leq j\leq n}\right)  ^{T}\right)  }_{=\left(  a_{j,n+1-i}\right)
_{1\leq i\leq n,\ 1\leq j\leq n}}=\det\left(  \left(  a_{j,n+1-i}\right)
_{1\leq i\leq n,\ 1\leq j\leq n}\right)  ,
\]
so that%
\begin{align*}
\det\left(  \left(  a_{j,n+1-i}\right)  _{1\leq i\leq n,\ 1\leq j\leq
n}\right)   &  =\det\left(  \left(  \left(  a_{i,n+1-j}\right)  _{1\leq i\leq
n,\ 1\leq j\leq n}\right)  ^{T}\right) \\
&  =\det\left(  \left(  a_{i,n+1-j}\right)  _{1\leq i\leq n,\ 1\leq j\leq
n}\right)
\end{align*}
(by Exercise \ref{exe.ps4.4}, applied to $A=\left(  a_{i,n+1-j}\right)
_{1\leq i\leq n,\ 1\leq j\leq n}$). Therefore,%
\begin{align*}
&  \det\left(  \left(  a_{i,n+1-j}\right)  _{1\leq i\leq n,\ 1\leq j\leq
n}\right) \\
&  =\det\left(  \left(  a_{j,n+1-i}\right)  _{1\leq i\leq n,\ 1\leq j\leq
n}\right)  =\left(  -1\right)  ^{n\left(  n-1\right)  /2}\underbrace{\det
\left(  \left(  a_{j,i}\right)  _{1\leq i\leq n,\ 1\leq j\leq n}\right)
}_{\substack{=\det\left(  \left(  a_{i,j}\right)  _{1\leq i\leq n,\ 1\leq
j\leq n}\right)  \\\text{(by (\ref{pf.lem.vander-det.lem-rearr-col.2}))}%
}}\ \ \ \ \ \ \ \ \ \ \left(  \text{by
(\ref{pf.lem.vander-det.lem-rearr-col.1})}\right) \\
&  =\left(  -1\right)  ^{n\left(  n-1\right)  /2}\det\left(  \left(
a_{i,j}\right)  _{1\leq i\leq n,\ 1\leq j\leq n}\right)  .
\end{align*}
This proves Lemma \ref{lem.vander-det.lem-rearr-col}.
\end{proof}

Our next lemma is a mild variation of Exercise \ref{exe.det.vdm-pol}:

\begin{lemma}
\label{lem.sol.det.kratt-lem6.vdm-pol-rev}Let $n\in\mathbb{N}$. Let
$a_{1},a_{2},\ldots,a_{n}\in\mathbb{K}$. Furthermore, for each $j\in\left\{
1,2,\ldots,n\right\}  $, let $Q_{j}\in\mathbb{K}\left[  X\right]  $ be a
polynomial such that%
\begin{equation}
\deg\left(  Q_{j}\right)  \leq n-j.
\label{eq.lem.sol.det.kratt-lem6.vdm-pol-rev.ass}%
\end{equation}
Then,%
\[
\det\left(  \left(  Q_{j}\left(  a_{i}\right)  \right)  _{1\leq i\leq
n,\ 1\leq j\leq n}\right)  =\left(  \prod_{j=1}^{n}\left[  X^{n-j}\right]
\left(  Q_{j}\right)  \right)  \cdot\prod_{1\leq i<j\leq n}\left(  a_{i}%
-a_{j}\right)  .
\]

\end{lemma}

\begin{proof}
[Proof of Lemma \ref{lem.sol.det.kratt-lem6.vdm-pol-rev}.]For each
$j\in\left\{  1,2,\ldots,n\right\}  $, we define a polynomial $P_{j}%
\in\mathbb{K}\left[  X\right]  $ by
\begin{equation}
P_{j}=Q_{n+1-j}. \label{pf.lem.sol.det.kratt-lem6.vdm-pol-rev.Pj=}%
\end{equation}
Then, for each $j\in\left\{  1,2,\ldots,n\right\}  $, the polynomial $P_{j}$
satisfies%
\[
\deg\left(  P_{j}\right)  \leq j-1
\]
\footnote{\textit{Proof.} Let $j\in\left\{  1,2,\ldots,n\right\}  $. Then,
$n+1-j\in\left\{  1,2,\ldots,n\right\}  $. However,
(\ref{pf.lem.sol.det.kratt-lem6.vdm-pol-rev.Pj=}) yields $P_{j}=Q_{n+1-j}$.
Therefore,%
\begin{align*}
\deg\left(  P_{j}\right)   &  =\deg\left(  Q_{n+1-j}\right) \\
&  \leq n-\left(  n+1-j\right)  \ \ \ \ \ \ \ \ \ \ \ \left(  \text{by
(\ref{eq.lem.sol.det.kratt-lem6.vdm-pol-rev.ass}), applied to }n+1-j\text{
instead of }j\right) \\
&  =j-1,
\end{align*}
qed.}. Hence, Lemma \ref{lem.sol.det.vdm-pol.1} yields%
\[
\det\left(  \left(  P_{j}\left(  a_{i}\right)  \right)  _{1\leq i\leq
n,\ 1\leq j\leq n}\right)  =\left(  \prod_{j=1}^{n}\left[  X^{j-1}\right]
\left(  P_{j}\right)  \right)  \cdot\det\left(  \left(  a_{i}^{j-1}\right)
_{1\leq i\leq n,\ 1\leq j\leq n}\right)  .
\]

However,
\begin{align*}
\det\left(  \left(  \underbrace{P_{j}}_{\substack{=Q_{n+1-j}\\\text{(by
(\ref{pf.lem.sol.det.kratt-lem6.vdm-pol-rev.Pj=}))}}}\left(  a_{i}\right)
\right)  _{1\leq i\leq n,\ 1\leq j\leq n}\right)   &  =\det\left(  \left(
Q_{n+1-j}\left(  a_{i}\right)  \right)  _{1\leq i\leq n,\ 1\leq j\leq
n}\right) \\
&  =\left(  -1\right)  ^{n\left(  n-1\right)  /2}\det\left(  \left(
Q_{j}\left(  a_{i}\right)  \right)  _{1\leq i\leq n,\ 1\leq j\leq n}\right)
\end{align*}
(by Lemma \ref{lem.vander-det.lem-rearr-col}, applied to $Q_{j}\left(
a_{i}\right)  $ instead of $a_{i,j}$). Hence,%
\begin{align*}
&  \left(  -1\right)  ^{n\left(  n-1\right)  /2}\det\left(  \left(
Q_{j}\left(  a_{i}\right)  \right)  _{1\leq i\leq n,\ 1\leq j\leq n}\right) \\
&  =\det\left(  \left(  P_{j}\left(  a_{i}\right)  \right)  _{1\leq i\leq
n,\ 1\leq j\leq n}\right) \\
&  =\left(  \prod_{j=1}^{n}\left[  X^{j-1}\right]  \left(  \underbrace{P_{j}%
}_{\substack{=Q_{n+1-j}\\\text{(by
(\ref{pf.lem.sol.det.kratt-lem6.vdm-pol-rev.Pj=}))}}}\right)  \right)
\cdot\det\left(  \left(  \underbrace{a_{i}^{j-1}}_{\substack{=a_{i}^{n-\left(
n+1-j\right)  }\\\text{(since }j-1=n-\left(  n+1-j\right)  \text{)}}}\right)
_{1\leq i\leq n,\ 1\leq j\leq n}\right) \\
&  =\underbrace{\left(  \prod_{j=1}^{n}\left[  X^{j-1}\right]  \left(
Q_{n+1-j}\right)  \right)  }_{\substack{=\prod_{j=1}^{n}\left[  X^{\left(
n+1-j\right)  -1}\right]  \left(  Q_{n+1-\left(  n+1-j\right)  }\right)
\\\text{(here, we have substituted }n+1-j\text{ for }j\\\text{in the
product)}}}\cdot\underbrace{\det\left(  \left(  a_{i}^{n-\left(  n+1-j\right)
}\right)  _{1\leq i\leq n,\ 1\leq j\leq n}\right)  }_{\substack{=\left(
-1\right)  ^{n\left(  n-1\right)  /2}\det\left(  \left(  a_{i}^{n-j}\right)
_{1\leq i\leq n,\ 1\leq j\leq n}\right)  \\\text{(by Lemma
\ref{lem.vander-det.lem-rearr-col},}\\\text{applied to }a_{i}^{n-j}\text{
instead of }a_{i,j}\text{)}}}\\
&  =\left(  \prod_{j=1}^{n}\underbrace{\left[  X^{\left(  n+1-j\right)
-1}\right]  \left(  Q_{n+1-\left(  n+1-j\right)  }\right)  }%
_{\substack{=\left[  X^{n-j}\right]  \left(  Q_{j}\right)  \\\text{(since
}\left(  n+1-j\right)  -1=n-j\\\text{and }n+1-\left(  n+1-j\right)
=j\text{)}}}\right)  \cdot\left(  -1\right)  ^{n\left(  n-1\right)
/2}\underbrace{\det\left(  \left(  a_{i}^{n-j}\right)  _{1\leq i\leq n,\ 1\leq
j\leq n}\right)  }_{\substack{=\prod_{1\leq i<j\leq n}\left(  a_{i}%
-a_{j}\right)  \\\text{(by Theorem \ref{thm.vander-det} \textbf{(a)}%
,}\\\text{applied to }x_{i}=a_{i}\text{)}}}\\
&  =\left(  \prod_{j=1}^{n}\left[  X^{n-j}\right]  \left(  Q_{j}\right)
\right)  \cdot\left(  -1\right)  ^{n\left(  n-1\right)  /2}\cdot\prod_{1\leq
i<j\leq n}\left(  a_{i}-a_{j}\right) \\
&  =\left(  -1\right)  ^{n\left(  n-1\right)  /2}\cdot\left(  \prod_{j=1}%
^{n}\left[  X^{n-j}\right]  \left(  Q_{j}\right)  \right)  \cdot\prod_{1\leq
i<j\leq n}\left(  a_{i}-a_{j}\right)  .
\end{align*}
Multiplying both sides of this equality by $\left(  -1\right)  ^{n\left(
n-1\right)  /2}$, we obtain%
\begin{align*}
&  \left(  -1\right)  ^{n\left(  n-1\right)  /2}\left(  -1\right)  ^{n\left(
n-1\right)  /2}\det\left(  \left(  Q_{j}\left(  a_{i}\right)  \right)  _{1\leq
i\leq n,\ 1\leq j\leq n}\right) \\
&  =\underbrace{\left(  -1\right)  ^{n\left(  n-1\right)  /2}\left(
-1\right)  ^{n\left(  n-1\right)  /2}}_{\substack{=\left(  \left(  -1\right)
^{n\left(  n-1\right)  /2}\right)  ^{2}=1\\\text{(since }\left(  -1\right)
^{n\left(  n-1\right)  /2}\text{ is either }-1\text{ or }1\text{)}}%
}\cdot\left(  \prod_{j=1}^{n}\left[  X^{n-j}\right]  \left(  Q_{j}\right)
\right)  \cdot\prod_{1\leq i<j\leq n}\left(  a_{i}-a_{j}\right) \\
&  =\left(  \prod_{j=1}^{n}\left[  X^{n-j}\right]  \left(  Q_{j}\right)
\right)  \cdot\prod_{1\leq i<j\leq n}\left(  a_{i}-a_{j}\right)  ,
\end{align*}
so that%
\begin{align*}
&  \left(  \prod_{j=1}^{n}\left[  X^{n-j}\right]  \left(  Q_{j}\right)
\right)  \cdot\prod_{1\leq i<j\leq n}\left(  a_{i}-a_{j}\right) \\
&  =\underbrace{\left(  -1\right)  ^{n\left(  n-1\right)  /2}\left(
-1\right)  ^{n\left(  n-1\right)  /2}}_{\substack{=\left(  \left(  -1\right)
^{n\left(  n-1\right)  /2}\right)  ^{2}=1\\\text{(since }\left(  -1\right)
^{n\left(  n-1\right)  /2}\text{ is either }-1\text{ or }1\text{)}}%
}\det\left(  \left(  Q_{j}\left(  a_{i}\right)  \right)  _{1\leq i\leq
n,\ 1\leq j\leq n}\right) \\
&  =\det\left(  \left(  Q_{j}\left(  a_{i}\right)  \right)  _{1\leq i\leq
n,\ 1\leq j\leq n}\right)  .
\end{align*}
This proves Lemma \ref{lem.sol.det.kratt-lem6.vdm-pol-rev}.
\end{proof}

We are now ready to solve Exercise \ref{exe.det.kratt-lem6}:

\begin{proof}
[Solution to Exercise \ref{exe.det.kratt-lem6}.]Let $j\in\left\{
1,2,\ldots,n\right\}  $. Then, $P_{j}$ is a polynomial such that $\deg\left(
P_{j}\right)  \leq j-1$ (by the definition of $P_{j}$).

For each $k\in\left\{  1,2,\ldots,j\right\}  $, let $R_{j,k}\in\mathbb{K}%
\left[  X\right]  $ be the polynomial%
\[
\left(  X-b_{k+1}\right)  \left(  X-b_{k+2}\right)  \cdots\left(
X-b_{j}\right)  .
\]
This polynomial $R_{j,k}$ is a product of $j-k$ monic polynomials of degree
$1$ (namely, of the $j-k$ polynomials $X-b_{k+1},X-b_{k+2},\ldots,X-b_{j}$),
and thus itself is a monic polynomial of degree $\underbrace{1+1+\cdots
+1}_{j-k\text{ times}}$ (since a product of several monic polynomials is
always a monic polynomial, and its degree is the sum of their degrees). In
other words, $R_{j,k}$ is a monic polynomial of degree $j-k$ (since
$\underbrace{1+1+\cdots+1}_{j-k\text{ times}}=\left(  j-k\right)  \cdot1=j-k$).

Hence, we can apply Lemma \ref{lem.sol.det.kratt-lem6.pol-div.1} to $R_{j,k}$
and $P_{j}$ instead of $R_{k}$ and $P$. We thus conclude that there exist $j$
elements $d_{1},d_{2},\ldots,d_{j}$ of $\mathbb{K}$ such that%
\[
P_{j}=d_{1}R_{j,1}+d_{2}R_{j,2}+\cdots+d_{j}R_{j,j}.
\]
Consider these $j$ elements $d_{1},d_{2},\ldots,d_{j}$, and denote them by
$c_{j,1},c_{j,2},\ldots,c_{j,j}$. Thus, $c_{j,1},c_{j,2},\ldots,c_{j,j}$ are
$j$ elements of $\mathbb{K}$ that satisfy%
\begin{equation}
P_{j}=c_{j,1}R_{j,1}+c_{j,2}R_{j,2}+\cdots+c_{j,j}R_{j,j}.
\label{sol.det.kratt-lem6.Pj=sum}%
\end{equation}

We extend the $j$-tuple $\left(  c_{j,1},c_{j,2},\ldots,c_{j,j}\right)
\in\mathbb{K}^{j}$ to an $n$-tuple $\left(  c_{j,1},c_{j,2},\ldots
,c_{j,n}\right)  \in\mathbb{K}^{n}$ by setting%
\begin{equation}
c_{j,k}=0\ \ \ \ \ \ \ \ \ \ \text{for each }k\in\left\{  j+1,j+2,\ldots
,n\right\}  . \label{sol.det.kratt-lem6.cji=0}%
\end{equation}

Now, (\ref{sol.det.kratt-lem6.Pj=sum}) becomes%
\begin{equation}
P_{j}=c_{j,1}R_{j,1}+c_{j,2}R_{j,2}+\cdots+c_{j,j}R_{j,j}=\sum_{k=1}%
^{j}c_{j,k}R_{j,k}. \label{sol.det.kratt-lem6.Pj=sum2}%
\end{equation}

Substituting $b_{j}$ for $X$ on both sides of this equality, we obtain%
\begin{align}
P_{j}\left(  b_{j}\right)   &  =\left(  \sum_{k=1}^{j}c_{j,k}R_{j,k}\right)
\left(  b_{j}\right)  =\sum_{k=1}^{j}c_{j,k}R_{j,k}\left(  b_{j}\right)
\nonumber\\
&  =\sum_{k=1}^{j-1}c_{j,k}R_{j,k}\left(  b_{j}\right)  +c_{j,j}R_{j,j}\left(
b_{j}\right)  \label{sol.det.kratt-lem6.Pj=sum2aj}%
\end{align}
(here, we have split off the addend for $k=j$ from the sum). However, each
$k\in\left\{  1,2,\ldots,j-1\right\}  $ satisfies
\begin{equation}
R_{j,k}\left(  b_{j}\right)  =0 \label{sol.det.kratt-lem6.3a}%
\end{equation}
\footnote{\textit{Proof of (\ref{sol.det.kratt-lem6.3a}):} Let $k\in\left\{
1,2,\ldots,j-1\right\}  $. Then, $k\leq j-1$, so that $k+1\leq j$. However,
the definition of $R_{j,k}$ yields
\[
R_{j,k}=\left(  X-b_{k+1}\right)  \left(  X-b_{k+2}\right)  \cdots\left(
X-b_{j}\right)  =\prod_{i=k+1}^{j}\left(  X-b_{i}\right)  .
\]
Substituting $b_{j}$ for $X$ on both sides of this equality, we obtain%
\begin{align*}
R_{j,k}\left(  b_{j}\right)   &  =\prod_{i=k+1}^{j}\left(  b_{j}-b_{i}\right)
=\underbrace{\left(  b_{j}-b_{j}\right)  }_{=0}\cdot\prod_{i=k+1}^{j-1}\left(
b_{j}-b_{i}\right) \\
&  \ \ \ \ \ \ \ \ \ \ \ \ \ \ \ \ \ \ \ \ \left(
\begin{array}
[c]{c}%
\text{here, we have split off the factor}\\
\text{for }i=j\text{ from the product, since }k+1\leq j
\end{array}
\right) \\
&  =0.
\end{align*}
This proves (\ref{sol.det.kratt-lem6.3a}).}. Furthermore, $R_{j,j}\left(
b_{j}\right)  =1$\ \ \ \ \footnote{\textit{Proof.} The definition of $R_{j,k}$
yields
\[
R_{j,j}=\left(  X-b_{j+1}\right)  \left(  X-b_{j+2}\right)  \cdots\left(
X-b_{j}\right)  =\left(  \text{empty product}\right)  =1.
\]
Substituting $b_{j}$ for $X$ on both sides of this equality, we obtain
$R_{j,j}\left(  b_{j}\right)  =1$.}. Hence,
(\ref{sol.det.kratt-lem6.Pj=sum2aj}) becomes%
\begin{align}
P_{j}\left(  b_{j}\right)   &  =\sum_{k=1}^{j-1}c_{j,k}\underbrace{R_{j,k}%
\left(  b_{j}\right)  }_{\substack{=0\\\text{(by (\ref{sol.det.kratt-lem6.3a}%
))}}}+c_{j,j}\underbrace{R_{j,j}\left(  b_{j}\right)  }_{=1}=\underbrace{\sum
_{k=1}^{j-1}c_{j,k}0}_{=0}+c_{j,j}\nonumber\\
&  =c_{j,j}. \label{sol.det.kratt-lem6.Pjbj=cjj}%
\end{align}

For each $k\in\left\{  1,2,\ldots,j\right\}  $, we have%
\begin{align}
&  \underbrace{R_{j,k}}_{\substack{=\left(  X-b_{k+1}\right)  \left(
X-b_{k+2}\right)  \cdots\left(  X-b_{j}\right)  \\\text{(by the definition of
}R_{j,k}\text{)}}}\ \ \underbrace{Q_{j}}_{=\left(  X-b_{j+1}\right)  \left(
X-b_{j+2}\right)  \cdots\left(  X-b_{n}\right)  }\nonumber\\
&  =\left(  X-b_{k+1}\right)  \left(  X-b_{k+2}\right)  \cdots\left(
X-b_{j}\right)  \cdot\left(  X-b_{j+1}\right)  \left(  X-b_{j+2}\right)
\cdots\left(  X-b_{n}\right) \nonumber\\
&  =\left(  X-b_{k+1}\right)  \left(  X-b_{k+2}\right)  \cdots\left(
X-b_{n}\right)  =Q_{k} \label{sol.det.kratt-lem6.RjkQj=Qk}%
\end{align}
(since $Q_{k}$ was defined to be $\left(  X-b_{k+1}\right)  \left(
X-b_{k+2}\right)  \cdots\left(  X-b_{n}\right)  $).

Next, let us multiply both sides of the equality
(\ref{sol.det.kratt-lem6.Pj=sum2}) by $Q_{j}$. We thus obtain
\[
P_{j}Q_{j}=\left(  \sum_{k=1}^{j}c_{j,k}R_{j,k}\right)  Q_{j}=\sum_{k=1}%
^{j}c_{j,k}\underbrace{R_{j,k}Q_{j}}_{\substack{=Q_{k}\\\text{(by
(\ref{sol.det.kratt-lem6.RjkQj=Qk}))}}}=\sum_{k=1}^{j}c_{j,k}Q_{k}.
\]
Comparing this with%
\begin{align*}
\sum_{k=1}^{n}c_{j,k}Q_{k}  &  =\sum_{k=1}^{j}c_{j,k}Q_{k}+\sum_{k=j+1}%
^{n}\underbrace{c_{j,k}}_{\substack{=0\\\text{(by
(\ref{sol.det.kratt-lem6.cji=0}))}}}Q_{k}\\
&  \ \ \ \ \ \ \ \ \ \ \ \ \ \ \ \ \ \ \ \ \left(  \text{here, we have split
the sum at }k=j\text{, since }1\leq j\leq n\right) \\
&  =\sum_{k=1}^{j}c_{j,k}Q_{k}+\underbrace{\sum_{k=j+1}^{n}0Q_{k}}_{=0}%
=\sum_{k=1}^{j}c_{j,k}Q_{k},
\end{align*}
we obtain%
\begin{equation}
P_{j}Q_{j}=\sum_{k=1}^{n}c_{j,k}Q_{k}. \label{sol.det.kratt-lem6.PjQj=sum}%
\end{equation}

Now, forget that we fixed $j$. Thus, for each $j\in\left\{  1,2,\ldots
,n\right\}  $, we have defined $n$ elements $c_{j,1},c_{j,2},\ldots,c_{j,n}$
of $\mathbb{K}$ that satisfy the equalities (\ref{sol.det.kratt-lem6.cji=0})
and (\ref{sol.det.kratt-lem6.Pjbj=cjj}) and (\ref{sol.det.kratt-lem6.PjQj=sum}%
). In other words, we have defined an element $c_{j,i}\in\mathbb{K}$ for each
$j\in\left\{  1,2,\ldots,n\right\}  $ and each $i\in\left\{  1,2,\ldots
,n\right\}  $. (These are altogether $n\cdot n=n^{2}$ many elements.)

Define an $n\times n$-matrix $C\in\mathbb{K}^{n\times n}$ by
\[
C=\left(  c_{j,i}\right)  _{1\leq i\leq n,\ 1\leq j\leq n}.
\]
We have $c_{j,i}=0$ for every $\left(  i,j\right)  \in\left\{  1,2,\ldots
,n\right\}  ^{2}$ satisfying $i>j$\ \ \ \ \footnote{\textit{Proof.} Let
$\left(  i,j\right)  \in\left\{  1,2,\ldots,n\right\}  ^{2}$ satisfy $i>j$.
From $i>j$, we obtain $i\geq j+1$ (since $i$ and $j$ are integers). However,
$i\in\left\{  1,2,\ldots,n\right\}  $ (since $\left(  i,j\right)  \in\left\{
1,2,\ldots,n\right\}  ^{2}$) and thus $i\leq n$. Combined with $i\geq j+1$,
this yields $i\in\left\{  j+1,j+2,\ldots,n\right\}  $. Therefore,
(\ref{sol.det.kratt-lem6.cji=0}) (applied to $k=i$) yields $c_{j,i}=0$. Qed.}.
Hence, Lemma \ref{lem.det.uptriangular} (applied to $C$ and $c_{j,i}$ instead
of $A$ and $a_{i,j}$) yields
\begin{align}
\det C  &  =c_{1,1}c_{2,2}\cdots c_{n,n}\ \ \ \ \ \ \ \ \ \ \left(
\text{since }C=\left(  c_{j,i}\right)  _{1\leq i\leq n,\ 1\leq j\leq n}\right)
\nonumber\\
&  =\prod_{j=1}^{n}\underbrace{c_{j,j}}_{\substack{=P_{j}\left(  b_{j}\right)
\\\text{(by (\ref{sol.det.kratt-lem6.Pjbj=cjj}))}}}=\prod_{j=1}^{n}%
P_{j}\left(  b_{j}\right)  . \label{sol.det.kratt-lem6.detC=}%
\end{align}

On the other hand, let us define an $n\times n$-matrix $B\in\mathbb{K}%
^{n\times n}$ by%
\[
B=\left(  Q_{j}\left(  a_{i}\right)  \right)  _{1\leq i\leq n,\ 1\leq j\leq
n}.
\]
Thus,%
\begin{align*}
\det\underbrace{B}_{=\left(  Q_{j}\left(  a_{i}\right)  \right)  _{1\leq i\leq
n,\ 1\leq j\leq n}}  &  =\det\left(  \left(  Q_{j}\left(  a_{i}\right)
\right)  _{1\leq i\leq n,\ 1\leq j\leq n}\right) \\
&  =\left(  \prod_{j=1}^{n}\left[  X^{n-j}\right]  \left(  Q_{j}\right)
\right)  \cdot\prod_{1\leq i<j\leq n}\left(  a_{i}-a_{j}\right)
\end{align*}
(by Lemma \ref{lem.sol.det.kratt-lem6.vdm-pol-rev}).

However, each $j\in\left\{  1,2,\ldots,n\right\}  $ satisfies $\left[
X^{n-j}\right]  \left(  Q_{j}\right)  =1$\ \ \ \ \footnote{\textit{Proof.} Let
$j\in\left\{  1,2,\ldots,n\right\}  $. Thus, $Q_{j}=\left(  X-b_{j+1}\right)
\left(  X-b_{j+2}\right)  \cdots\left(  X-b_{n}\right)  $ (by the definition
of $Q_{j}$). Hence, the polynomial $Q_{j}$ is a product of $n-j$ monic
polynomials of degree $1$ (namely, of the $n-j$ polynomials $X-b_{j+1}%
,X-b_{j+2},\ldots,X-b_{n}$), and thus itself is a monic polynomial of degree
$\underbrace{1+1+\cdots+1}_{n-j\text{ times}}$ (since a product of several
monic polynomials is always a monic polynomial, and its degree is the sum of
their degrees). In other words, $Q_{j}$ is a monic polynomial of degree $n-j$
(since $\underbrace{1+1+\cdots+1}_{n-j\text{ times}}=\left(  n-j\right)
\cdot1=n-j$). Thus, the coefficient of $X^{n-j}$ in $Q_{j}$ is $1$. In other
words, $\left[  X^{n-j}\right]  \left(  Q_{j}\right)  =1$ (since $\left[
X^{n-j}\right]  \left(  Q_{j}\right)  =1$ was defined to be the coefficient of
$X^{n-j}$ in $Q_{j}$). Qed.}. Multiplying these equalities over all
$j\in\left\{  1,2,\ldots,n\right\}  $, we obtain $\prod_{j=1}^{n}\left[
X^{n-j}\right]  \left(  Q_{j}\right)  =\prod_{j=1}^{n}1=1$. Now, our above
computation becomes%
\begin{align}
\det B  &  =\underbrace{\left(  \prod_{j=1}^{n}\left[  X^{n-j}\right]  \left(
Q_{j}\right)  \right)  }_{=1}\cdot\prod_{1\leq i<j\leq n}\left(  a_{i}%
-a_{j}\right) \nonumber\\
&  =\prod_{1\leq i<j\leq n}\left(  a_{i}-a_{j}\right)  .
\label{sol.det.kratt-lem6.detB=}%
\end{align}

Now, we claim that%
\begin{equation}
\left(  P_{j}\left(  a_{i}\right)  Q_{j}\left(  a_{i}\right)  \right)  _{1\leq
i\leq n,\ 1\leq j\leq n}=BC. \label{sol.det.kratt-lem6.=BC}%
\end{equation}

[\textit{Proof of (\ref{sol.det.kratt-lem6.=BC}):} Let $\left(  i,j\right)
\in\left\{  1,2,\ldots,n\right\}  ^{2}$. Substituting $a_{i}$ for $X$ on both
sides of the equality (\ref{sol.det.kratt-lem6.PjQj=sum}), we obtain%
\[
\left(  P_{j}Q_{j}\right)  \left(  a_{i}\right)  =\left(  \sum_{k=1}%
^{n}c_{j,k}Q_{k}\right)  \left(  a_{i}\right)  =\sum_{k=1}^{n}%
\underbrace{c_{j,k}Q_{k}\left(  a_{i}\right)  }_{=Q_{k}\left(  a_{i}\right)
c_{j,k}}=\sum_{k=1}^{n}Q_{k}\left(  a_{i}\right)  c_{j,k}.
\]
In view of $\left(  P_{j}Q_{j}\right)  \left(  a_{i}\right)  =P_{j}\left(
a_{i}\right)  Q_{j}\left(  a_{i}\right)  $, this rewrites as%
\begin{equation}
P_{j}\left(  a_{i}\right)  Q_{j}\left(  a_{i}\right)  =\sum_{k=1}^{n}%
Q_{k}\left(  a_{i}\right)  c_{j,k}. \label{sol.det.kratt-lem6.=BC.pf.1}%
\end{equation}

Forget that we fixed $\left(  i,j\right)  $. We thus have proved that
(\ref{sol.det.kratt-lem6.=BC.pf.1}) holds for each $\left(  i,j\right)
\in\left\{  1,2,\ldots,n\right\}  ^{2}$. In other words,%
\[
\left(  P_{j}\left(  a_{i}\right)  Q_{j}\left(  a_{i}\right)  \right)  _{1\leq
i\leq n,\ 1\leq j\leq n}=\left(  \sum_{k=1}^{n}Q_{k}\left(  a_{i}\right)
c_{j,k}\right)  _{1\leq i\leq n,\ 1\leq j\leq n}.
\]
On the other hand, the definition of the product of two matrices yields%
\[
BC=\left(  \sum_{k=1}^{n}Q_{k}\left(  a_{i}\right)  c_{j,k}\right)  _{1\leq
i\leq n,\ 1\leq j\leq n}%
\]
(since $B=\left(  Q_{j}\left(  a_{i}\right)  \right)  _{1\leq i\leq n,\ 1\leq
j\leq n}$ and $C=\left(  c_{j,i}\right)  _{1\leq i\leq n,\ 1\leq j\leq n}$).
Comparing these two equalities, we find $\left(  P_{j}\left(  a_{i}\right)
Q_{j}\left(  a_{i}\right)  \right)  _{1\leq i\leq n,\ 1\leq j\leq n}=BC$. This
proves (\ref{sol.det.kratt-lem6.=BC}).]

Now,%
\begin{align*}
&  \det\left(  \underbrace{\left(  P_{j}\left(  a_{i}\right)  Q_{j}\left(
a_{i}\right)  \right)  _{1\leq i\leq n,\ 1\leq j\leq n}}%
_{\substack{=BC\\\text{(by (\ref{sol.det.kratt-lem6.=BC}))}}}\right) \\
&  =\det\left(  BC\right)  =\det B\cdot\det C\\
&  \ \ \ \ \ \ \ \ \ \ \left(  \text{by Theorem \ref{thm.det(AB)}, applied to
}B\text{ and }C\text{ instead of }A\text{ and }B\right) \\
&  =\underbrace{\det C}_{\substack{=\prod_{j=1}^{n}P_{j}\left(  b_{j}\right)
\\\text{(by (\ref{sol.det.kratt-lem6.detC=}))}}}\cdot\underbrace{\det
B}_{\substack{=\prod_{1\leq i<j\leq n}\left(  a_{i}-a_{j}\right)  \\\text{(by
(\ref{sol.det.kratt-lem6.detB=}))}}}\\
&  =\left(  \prod_{j=1}^{n}P_{j}\left(  b_{j}\right)  \right)  \cdot
\prod_{1\leq i<j\leq n}\left(  a_{i}-a_{j}\right)  .
\end{align*}
Thus, Exercise \ref{exe.det.kratt-lem6} is solved.
\end{proof}

\subsection{Solution to Exercise \ref{exe.det.kratt-lem6-app}}

\begin{proof}
[Solution to Exercise \ref{exe.det.kratt-lem6-app}.]\textbf{(a)} For each
$j\in\left\{  1,2,\ldots,n\right\}  $, we define a polynomial $Q_{j}%
\in\mathbb{K}\left[  X\right]  $ by setting%
\begin{equation}
Q_{j}=\left(  X-b_{j+1}\right)  \left(  X-b_{j+2}\right)  \cdots\left(
X-b_{n}\right)  . \label{sol.det.kratt-lem6-app.a.Qj=}%
\end{equation}
For each $j\in\left\{  1,2,\ldots,n\right\}  $, we define a polynomial
$P_{j}\in\mathbb{K}\left[  X\right]  $ by setting%
\begin{equation}
P_{j}=\left(  X-c_{1}\right)  \left(  X-c_{2}\right)  \cdots\left(
X-c_{j-1}\right)  . \label{sol.det.kratt-lem6-app.a.Pj=}%
\end{equation}
This polynomial $P_{j}$ satisfies $\deg\left(  P_{j}\right)  \leq
j-1$\ \ \ \ \footnote{\textit{Proof.} Let $j\in\left\{  1,2,\ldots,n\right\}
$. Then, $P_{j}=\left(  X-c_{1}\right)  \left(  X-c_{2}\right)  \cdots\left(
X-c_{j-1}\right)  $. Hence, this polynomial $P_{j}$ is a product of $j-1$
monic polynomials of degree $1$ (namely, of the $j-1$ polynomials
$X-c_{1},X-c_{2},\ldots,X-c_{j-1}$), and thus itself is a monic polynomial of
degree $\underbrace{1+1+\cdots+1}_{j-1\text{ times}}$ (since a product of
several monic polynomials is always a monic polynomial, and its degree is the
sum of their degrees). In other words, $P_{j}$ is a monic polynomial of degree
$j-1$ (since $\underbrace{1+1+\cdots+1}_{j-1\text{ times}}=\left(  j-1\right)
\cdot1=j-1$). Hence, in particular, this polynomial $P_{j}$ has degree $j-1$,
so that $\deg\left(  P_{j}\right)  =j-1\leq j-1$. Qed.}. Thus, Exercise
\ref{exe.det.kratt-lem6} yields
\begin{align}
&  \det\left(  \left(  P_{j}\left(  a_{i}\right)  Q_{j}\left(  a_{i}\right)
\right)  _{1\leq i\leq n,\ 1\leq j\leq n}\right) \nonumber\\
&  =\left(  \prod_{j=1}^{n}P_{j}\left(  b_{j}\right)  \right)  \cdot
\prod_{1\leq i<j\leq n}\left(  a_{i}-a_{j}\right)  .
\label{sol.det.kratt-lem6-app.a.1}%
\end{align}

However, each $\left(  i,j\right)  \in\left\{  1,2,\ldots,n\right\}  ^{2}$
satisfies%
\begin{align}
&  \underbrace{P_{j}\left(  a_{i}\right)  }_{\substack{=\left(  a_{i}%
-c_{1}\right)  \left(  a_{i}-c_{2}\right)  \cdots\left(  a_{i}-c_{j-1}\right)
\\\text{(here, we have substituted }a_{i}\text{ for }X\\\text{on both sides of
the equality (\ref{sol.det.kratt-lem6-app.a.Pj=}))}}}\ \ \underbrace{Q_{j}%
\left(  a_{i}\right)  }_{\substack{=\left(  a_{i}-b_{j+1}\right)  \left(
a_{i}-b_{j+2}\right)  \cdots\left(  a_{i}-b_{n}\right)  \\\text{(here, we have
substituted }a_{i}\text{ for }X\\\text{on both sides of the equality
(\ref{sol.det.kratt-lem6-app.a.Qj=}))}}}\nonumber\\
&  =\underbrace{\left(  a_{i}-c_{1}\right)  \left(  a_{i}-c_{2}\right)
\cdots\left(  a_{i}-c_{j-1}\right)  }_{=\prod_{u=1}^{j-1}\left(  a_{i}%
-c_{u}\right)  }\cdot\underbrace{\left(  a_{i}-b_{j+1}\right)  \left(
a_{i}-b_{j+2}\right)  \cdots\left(  a_{i}-b_{n}\right)  }_{=\prod_{u=j+1}%
^{n}\left(  a_{i}-b_{u}\right)  }\nonumber\\
&  =\left(  \prod_{u=1}^{j-1}\left(  a_{i}-c_{u}\right)  \right)  \left(
\prod_{u=j+1}^{n}\left(  a_{i}-b_{u}\right)  \right)  .
\label{sol.det.kratt-lem6-app.a.PjQjai}%
\end{align}

Furthermore, each $j\in\left\{  1,2,\ldots,n\right\}  $ satisfies%
\[
P_{j}\left(  b_{j}\right)  =\left(  b_{j}-c_{1}\right)  \left(  b_{j}%
-c_{2}\right)  \cdots\left(  b_{j}-c_{j-1}\right)
\]
(here, we have substituted $b_{j}$ for $X$ on both sides of the equality
(\ref{sol.det.kratt-lem6-app.a.Pj=})). Multiplying these equalities over all
$j\in\left\{  1,2,\ldots,n\right\}  $, we obtain%
\begin{align*}
\prod_{j=1}^{n}P_{j}\left(  b_{j}\right)   &  =\prod_{j=1}^{n}%
\underbrace{\left(  \left(  b_{j}-c_{1}\right)  \left(  b_{j}-c_{2}\right)
\cdots\left(  b_{j}-c_{j-1}\right)  \right)  }_{=\prod_{i=1}^{j-1}\left(
b_{j}-c_{i}\right)  }=\underbrace{\prod_{j=1}^{n}\ \ \prod_{i=1}^{j-1}%
}_{=\prod_{1\leq i<j\leq n}}\left(  b_{j}-c_{i}\right) \\
&  =\prod_{1\leq i<j\leq n}\left(  b_{j}-c_{i}\right)  .
\end{align*}
In view of this, we can rewrite (\ref{sol.det.kratt-lem6-app.a.1}) as%
\begin{align*}
&  \det\left(  \left(  P_{j}\left(  a_{i}\right)  Q_{j}\left(  a_{i}\right)
\right)  _{1\leq i\leq n,\ 1\leq j\leq n}\right) \\
&  =\left(  \prod_{1\leq i<j\leq n}\left(  b_{j}-c_{i}\right)  \right)
\cdot\prod_{1\leq i<j\leq n}\left(  a_{i}-a_{j}\right)  .
\end{align*}
In view of (\ref{sol.det.kratt-lem6-app.a.PjQjai}), this furthermore rewrites
as%
\begin{align*}
&  \det\left(  \left(  \left(  \prod_{u=1}^{j-1}\left(  a_{i}-c_{u}\right)
\right)  \left(  \prod_{u=j+1}^{n}\left(  a_{i}-b_{u}\right)  \right)
\right)  _{1\leq i\leq n,\ 1\leq j\leq n}\right) \\
&  =\left(  \prod_{1\leq i<j\leq n}\left(  b_{j}-c_{i}\right)  \right)
\cdot\prod_{1\leq i<j\leq n}\left(  a_{i}-a_{j}\right)  .
\end{align*}
This solves Exercise \ref{exe.det.kratt-lem6-app} \textbf{(a)}.

\textbf{(b)} Exercise \ref{exe.det.kratt-lem6-app} \textbf{(a)} (applied to
$-b_{j}$ and $-c_{i}$ instead of $b_{j}$ and $c_{i}$) yields%
\begin{align*}
&  \det\left(  \left(  \left(  \prod_{u=1}^{j-1}\left(  a_{i}-\left(
-c_{u}\right)  \right)  \right)  \left(  \prod_{u=j+1}^{n}\left(
a_{i}-\left(  -b_{u}\right)  \right)  \right)  \right)  _{1\leq i\leq
n,\ 1\leq j\leq n}\right) \\
&  =\left(  \prod_{1\leq i<j\leq n}\underbrace{\left(  \left(  -b_{j}\right)
-\left(  -c_{i}\right)  \right)  }_{=c_{i}-b_{j}}\right)  \cdot\prod_{1\leq
i<j\leq n}\left(  a_{i}-a_{j}\right) \\
&  =\left(  \prod_{1\leq i<j\leq n}\left(  c_{i}-b_{j}\right)  \right)
\cdot\prod_{1\leq i<j\leq n}\left(  a_{i}-a_{j}\right)  .
\end{align*}
Comparing this with%
\begin{align*}
&  \det\left(  \left(  \left(  \prod_{u=1}^{j-1}\underbrace{\left(
a_{i}-\left(  -c_{u}\right)  \right)  }_{=a_{i}+c_{u}}\right)  \left(
\prod_{u=j+1}^{n}\underbrace{\left(  a_{i}-\left(  -b_{u}\right)  \right)
}_{=a_{i}+b_{u}}\right)  \right)  _{1\leq i\leq n,\ 1\leq j\leq n}\right) \\
&  =\det\left(  \left(  \left(  \prod_{u=1}^{j-1}\left(  a_{i}+c_{u}\right)
\right)  \left(  \prod_{u=j+1}^{n}\left(  a_{i}+b_{u}\right)  \right)
\right)  _{1\leq i\leq n,\ 1\leq j\leq n}\right)  ,
\end{align*}
we obtain%
\begin{align*}
&  \det\left(  \left(  \left(  \prod_{u=1}^{j-1}\left(  a_{i}+c_{u}\right)
\right)  \left(  \prod_{u=j+1}^{n}\left(  a_{i}+b_{u}\right)  \right)
\right)  _{1\leq i\leq n,\ 1\leq j\leq n}\right) \\
&  =\left(  \prod_{1\leq i<j\leq n}\left(  c_{i}-b_{j}\right)  \right)
\cdot\prod_{1\leq i<j\leq n}\left(  a_{i}-a_{j}\right)  .
\end{align*}
This solves Exercise \ref{exe.det.kratt-lem6-app} \textbf{(b)}.

\textbf{(c)} Here is a second solution to Exercise \ref{exe.cauchy-det}:

\begin{proof}
[Second solution to Exercise \ref{exe.cauchy-det}.]Exercise
\ref{exe.det.kratt-lem6-app} \textbf{(b)} (applied to $a_{i}=x_{i}$ and
$b_{j}=y_{j}$ and $c_{j}=y_{j}$) yields%
\begin{align}
&  \det\left(  \left(  \left(  \prod_{u=1}^{j-1}\left(  x_{i}+y_{u}\right)
\right)  \left(  \prod_{u=j+1}^{n}\left(  x_{i}+y_{u}\right)  \right)
\right)  _{1\leq i\leq n,\ 1\leq j\leq n}\right) \nonumber\\
&  =\left(  \prod_{1\leq i<j\leq n}\left(  y_{i}-y_{j}\right)  \right)
\cdot\prod_{1\leq i<j\leq n}\left(  x_{i}-x_{j}\right) \nonumber\\
&  =\prod_{1\leq i<j\leq n}\underbrace{\left(  \left(  y_{i}-y_{j}\right)
\left(  x_{i}-x_{j}\right)  \right)  }_{=\left(  x_{i}-x_{j}\right)  \left(
y_{i}-y_{j}\right)  }\nonumber\\
&  =\prod_{1\leq i<j\leq n}\left(  \left(  x_{i}-x_{j}\right)  \left(
y_{i}-y_{j}\right)  \right)  . \label{sol.det.kratt-lem-6-app.c.1}%
\end{align}

Now, for each $i\in\left\{  1,2,\ldots,n\right\}  $, set%
\[
d_{i}=\prod_{u=1}^{n}\left(  x_{i}+y_{u}\right)  .
\]
Furthermore, for each $\left(  i,j\right)  \in\left\{  1,2,\ldots,n\right\}
^{2}$, we set
\[
a_{i,j}=\dfrac{d_{i}}{x_{i}+y_{j}}.
\]

Let $\left(  i,j\right)  \in\left\{  1,2,\ldots,n\right\}  ^{2}$. Then,%
\begin{align}
d_{i}  &  =\prod_{u=1}^{n}\left(  x_{i}+y_{u}\right)  =\underbrace{\left(
\prod_{u=1}^{j}\left(  x_{i}+y_{u}\right)  \right)  }_{\substack{=\left(
x_{i}+y_{j}\right)  \cdot\prod_{u=1}^{j-1}\left(  x_{i}+y_{u}\right)
\\\text{(here, we have split off the factor}\\\text{for }u=j\text{ from the
product)}}}\cdot\left(  \prod_{u=j+1}^{n}\left(  x_{i}+y_{u}\right)  \right)
\nonumber\\
&  \ \ \ \ \ \ \ \ \ \ \ \ \ \ \ \ \ \ \ \ \left(  \text{here, we have split
the product at }u=j\right) \nonumber\\
&  =\left(  x_{i}+y_{j}\right)  \cdot\left(  \prod_{u=1}^{j-1}\left(
x_{i}+y_{u}\right)  \right)  \cdot\left(  \prod_{u=j+1}^{n}\left(  x_{i}%
+y_{u}\right)  \right)  . \label{sol.det.kratt-lem6-app.c.3}%
\end{align}
Now, the definition of $a_{i,j}$ yields%
\begin{equation}
a_{i,j}=\dfrac{d_{i}}{x_{i}+y_{j}}=\left(  \prod_{u=1}^{j-1}\left(
x_{i}+y_{u}\right)  \right)  \cdot\left(  \prod_{u=j+1}^{n}\left(  x_{i}%
+y_{u}\right)  \right)  \label{sol.det.kratt-lem6-app.c.4}%
\end{equation}
(here, we have divided both sides of the equality
(\ref{sol.det.kratt-lem6-app.c.3}) by $x_{i}+y_{j}$).

Forget that we fixed $\left(  i,j\right)  $. We thus have proved
(\ref{sol.det.kratt-lem6-app.c.4}) for each $\left(  i,j\right)  \in\left\{
1,2,\ldots,n\right\}  ^{2}$. In other words,%
\[
\left(  a_{i,j}\right)  _{1\leq i\leq n,\ 1\leq j\leq n}=\left(  \left(
\prod_{u=1}^{j-1}\left(  x_{i}+y_{u}\right)  \right)  \cdot\left(
\prod_{u=j+1}^{n}\left(  x_{i}+y_{u}\right)  \right)  \right)  _{1\leq i\leq
n,\ 1\leq j\leq n}.
\]
Hence,%
\begin{align*}
&  \det\left(  \left(  a_{i,j}\right)  _{1\leq i\leq n,\ 1\leq j\leq n}\right)
\\
&  =\det\left(  \left(  \left(  \prod_{u=1}^{j-1}\left(  x_{i}+y_{u}\right)
\right)  \cdot\left(  \prod_{u=j+1}^{n}\left(  x_{i}+y_{u}\right)  \right)
\right)  _{1\leq i\leq n,\ 1\leq j\leq n}\right) \\
&  =\prod_{1\leq i<j\leq n}\left(  \left(  x_{i}-x_{j}\right)  \left(
y_{i}-y_{j}\right)  \right)  \ \ \ \ \ \ \ \ \ \ \left(  \text{by
(\ref{sol.det.kratt-lem-6-app.c.1})}\right)  .
\end{align*}

Now, Lemma \ref{lem.sol.cauchy-det-lem.1} (applied to $c_{i}=\dfrac{1}{d_{i}}%
$) yields%
\begin{align*}
\det\left(  \left(  \dfrac{1}{d_{i}}a_{i,j}\right)  _{1\leq i\leq n,\ 1\leq
j\leq n}\right)   &  =\underbrace{\left(  \prod_{i=1}^{n}\dfrac{1}{d_{i}%
}\right)  }_{=\dfrac{1}{\prod_{i=1}^{n}d_{i}}}\cdot\underbrace{\det\left(
\left(  a_{i,j}\right)  _{1\leq i\leq n,\ 1\leq j\leq n}\right)  }%
_{=\prod_{1\leq i<j\leq n}\left(  \left(  x_{i}-x_{j}\right)  \left(
y_{i}-y_{j}\right)  \right)  }\\
&  =\dfrac{1}{\prod_{i=1}^{n}d_{i}}\cdot\prod_{1\leq i<j\leq n}\left(  \left(
x_{i}-x_{j}\right)  \left(  y_{i}-y_{j}\right)  \right) \\
&  =\dfrac{\prod_{1\leq i<j\leq n}\left(  \left(  x_{i}-x_{j}\right)  \left(
y_{i}-y_{j}\right)  \right)  }{\prod_{i=1}^{n}d_{i}}.
\end{align*}
In view of%
\begin{align*}
\prod_{i=1}^{n}\underbrace{d_{i}}_{\substack{=\prod_{u=1}^{n}\left(
x_{i}+y_{u}\right)  \\\text{(by the definition of }d_{i}\text{)}}}  &
=\prod_{i=1}^{n}\ \ \prod_{u=1}^{n}\left(  x_{i}+y_{u}\right)
=\underbrace{\prod_{i=1}^{n}\ \ \prod_{j=1}^{n}}_{=\prod_{\left(  i,j\right)
\in\left\{  1,2,\ldots,n\right\}  ^{2}}}\left(  x_{i}+y_{j}\right) \\
&  \ \ \ \ \ \ \ \ \ \ \ \ \ \ \ \ \ \ \ \ \left(
\begin{array}
[c]{c}%
\text{here, we have renamed the index }u\\
\text{as }j\text{ in the second product}%
\end{array}
\right) \\
&  =\prod_{\left(  i,j\right)  \in\left\{  1,2,\ldots,n\right\}  ^{2}}\left(
x_{i}+y_{j}\right)  ,
\end{align*}
we can rewrite this as
\begin{equation}
\det\left(  \left(  \dfrac{1}{d_{i}}a_{i,j}\right)  _{1\leq i\leq n,\ 1\leq
j\leq n}\right)  =\dfrac{\prod_{1\leq i<j\leq n}\left(  \left(  x_{i}%
-x_{j}\right)  \left(  y_{i}-y_{j}\right)  \right)  }{\prod_{\left(
i,j\right)  \in\left\{  1,2,\ldots,n\right\}  ^{2}}\left(  x_{i}+y_{j}\right)
}. \label{sol.det.kratt-lem-6-app.c.9}%
\end{equation}

However, each $\left(  i,j\right)  \in\left\{  1,2,\ldots,n\right\}  ^{2}$
satisfies%
\[
\dfrac{1}{d_{i}}\underbrace{a_{i,j}}_{\substack{=\dfrac{d_{i}}{x_{i}+y_{j}%
}\\\text{(by the definition of }a_{i,j}\text{)}}}=\dfrac{1}{d_{i}}\cdot
\dfrac{d_{i}}{x_{i}+y_{j}}=\dfrac{1}{x_{i}+y_{j}}.
\]
Thus, we can rewrite (\ref{sol.det.kratt-lem-6-app.c.9}) as%
\[
\det\left(  \left(  \dfrac{1}{x_{i}+y_{j}}\right)  _{1\leq i\leq n,\ 1\leq
j\leq n}\right)  =\dfrac{\prod_{1\leq i<j\leq n}\left(  \left(  x_{i}%
-x_{j}\right)  \left(  y_{i}-y_{j}\right)  \right)  }{\prod_{\left(
i,j\right)  \in\left\{  1,2,\ldots,n\right\}  ^{2}}\left(  x_{i}+y_{j}\right)
}.
\]
Thus, Exercise \ref{exe.cauchy-det} is solved again.
\end{proof}

This solves Exercise \ref{exe.det.kratt-lem6-app} \textbf{(c)}.
\end{proof}

\subsection{Solution to Exercise \ref{exe.det.fraser-yeats}}

Before we solve Exercise \ref{exe.det.fraser-yeats}, we shall prove a lemma:

\begin{lemma}
\label{lem.sol.det.fraser-yeats.1}Let $n\in\mathbb{N}$. Let $A=\left(
a_{i,j}\right)  _{1\leq i\leq n,\ 1\leq j\leq n}\in\mathbb{K}^{n\times n}$ be
an $n\times n$-matrix. Let $k\in\mathbb{N}$. Let $u_{1},u_{2},\ldots,u_{k}$ be
$k$ elements of $\left\{  1,2,\ldots,n\right\}  $. Let $v_{1},v_{2}%
,\ldots,v_{k}$ be $k$ further elements of $\left\{  1,2,\ldots,n\right\}  $.
Let $s\in\left\{  1,2,\ldots,n\right\}  $ and $t\in\left\{  1,2,\ldots
,n\right\}  $. Then,%
\begin{align*}
&  \det\left(  \operatorname*{sub}\nolimits_{s,u_{1},u_{2},\ldots,u_{k}%
}^{t,v_{1},v_{2},\ldots,v_{k}}A\right) \\
&  =a_{s,t}\det\left(  \operatorname*{sub}\nolimits_{u_{1},u_{2},\ldots,u_{k}%
}^{v_{1},v_{2},\ldots,v_{k}}A\right)  +\sum_{p=1}^{k}\left(  -1\right)
^{p}a_{s,v_{p}}\det\left(  \operatorname*{sub}\nolimits_{u_{1},u_{2}%
,\ldots,u_{k}}^{t,v_{1},v_{2},\ldots,\widehat{v_{p}},\ldots,v_{k}}A\right)  .
\end{align*}
Here, we are using the notation introduced in Definition \ref{def.hat-omit}
(so that \textquotedblleft$t,v_{1},v_{2},\ldots,\widehat{v_{p}},\ldots,v_{k}%
$\textquotedblright\ means \textquotedblleft$t,v_{1},v_{2},\ldots
,v_{p-1},v_{p+1},v_{p+2},\ldots,v_{k}$\textquotedblright), as well as the
notation of Definition \ref{def.submatrix}.
\end{lemma}

\begin{proof}
[Proof of Lemma \ref{lem.sol.det.fraser-yeats.1}.]Essentially, this follows by
expanding the determinant of the matrix $\operatorname*{sub}\nolimits_{s,u_{1}%
,u_{2},\ldots,u_{k}}^{t,v_{1},v_{2},\ldots,v_{k}}A$ along the first row (using
Theorem \ref{thm.laplace.gen} \textbf{(a)}). Here are the details:

Set%
\[
B=\operatorname*{sub}\nolimits_{s,u_{1},u_{2},\ldots,u_{k}}^{t,v_{1}%
,v_{2},\ldots,v_{k}}A.
\]
Thus, $B$ is a $\left(  k+1\right)  \times\left(  k+1\right)  $-matrix.

We extend the $k$-tuple $\left(  u_{1},u_{2},\ldots,u_{k}\right)  $ to a
$\left(  k+1\right)  $-tuple $\left(  u_{0},u_{1},\ldots,u_{k}\right)  $ by
setting $u_{0}=s$. We extend the $k$-tuple $\left(  v_{1},v_{2},\ldots
,v_{k}\right)  $ to a $\left(  k+1\right)  $-tuple $\left(  v_{0},v_{1}%
,\ldots,v_{k}\right)  $ by setting $v_{0}=t$. We have%
\[
\left(  u_{0},u_{1},\ldots,u_{k}\right)  =\left(  \underbrace{u_{0}}%
_{=s},u_{1},u_{2},\ldots,u_{k}\right)  =\left(  s,u_{1},u_{2},\ldots
,u_{k}\right)
\]
and%
\[
\left(  v_{0},v_{1},\ldots,v_{k}\right)  =\left(  \underbrace{v_{0}}%
_{=t},v_{1},v_{2},\ldots,v_{k}\right)  =\left(  t,v_{1},v_{2},\ldots
,v_{k}\right)  .
\]
Therefore,%
\[
\operatorname*{sub}\nolimits_{u_{0},u_{1},\ldots,u_{k}}^{v_{0},v_{1}%
,\ldots,v_{k}}A=\operatorname*{sub}\nolimits_{s,u_{1},u_{2},\ldots,u_{k}%
}^{t,v_{1},v_{2},\ldots,v_{k}}A=B.
\]
Hence,%
\begin{align*}
B  &  =\operatorname*{sub}\nolimits_{u_{0},u_{1},\ldots,u_{k}}^{v_{0}%
,v_{1},\ldots,v_{k}}A\\
&  =\left(  a_{u_{x-1},v_{y-1}}\right)  _{1\leq x\leq k+1,\ 1\leq y\leq
k+1}\ \ \ \ \ \ \ \ \ \ \left(  \text{by Definition \ref{def.submatrix}, since
}A=\left(  a_{i,j}\right)  _{1\leq i\leq n,\ 1\leq j\leq n}\right) \\
&  =\left(  a_{u_{i-1},v_{j-1}}\right)  _{1\leq i\leq k+1,\ 1\leq j\leq k+1}%
\end{align*}
(here, we have renamed the index $\left(  x,y\right)  $ as $\left(
i,j\right)  $). Hence, Theorem \ref{thm.laplace.gen} \textbf{(a)} (applied to
$k+1$, $B$, $a_{u_{i-1},v_{j-1}}$ and $1$ instead of $n$, $A$, $a_{i,j}$ and
$p$) yields%
\begin{align}
\det B  &  =\sum_{q=1}^{k+1}\left(  -1\right)  ^{1+q}a_{u_{1-1},v_{q-1}}%
\det\left(  B_{\sim1,\sim q}\right) \nonumber\\
&  =\sum_{q=1}^{k+1}\left(  -1\right)  ^{1+q}a_{s,v_{q-1}}\det\left(
B_{\sim1,\sim q}\right)  \ \ \ \ \ \ \ \ \ \ \left(  \text{since }%
u_{1-1}=u_{0}=s\right) \nonumber\\
&  =\underbrace{\left(  -1\right)  ^{1+1}}_{=\left(  -1\right)  ^{2}%
=1}\underbrace{a_{s,v_{1-1}}}_{\substack{=a_{s,t}\\\text{(since }v_{1-1}%
=v_{0}=t\text{)}}}\det\left(  B_{\sim1,\sim1}\right)  +\sum_{q=2}^{k+1}\left(
-1\right)  ^{1+q}a_{s,v_{q-1}}\det\left(  B_{\sim1,\sim q}\right) \nonumber\\
&  \ \ \ \ \ \ \ \ \ \ \ \ \ \ \ \ \ \ \ \ \left(  \text{here, we have split
off the addend for }q=1\text{ from the sum}\right) \nonumber\\
&  =a_{s,t}\det\left(  B_{\sim1,\sim1}\right)  +\sum_{q=2}^{k+1}\left(
-1\right)  ^{1+q}a_{s,v_{q-1}}\det\left(  B_{\sim1,\sim q}\right) \nonumber\\
&  =a_{s,t}\det\left(  B_{\sim1,\sim1}\right)  +\sum_{p=1}^{k}%
\underbrace{\left(  -1\right)  ^{1+\left(  p+1\right)  }}_{\substack{=\left(
-1\right)  ^{p}\\\text{(since }1+\left(  p+1\right)  =p+2\equiv
p\operatorname{mod}2\text{)}}}\underbrace{a_{s,v_{\left(  p+1\right)  -1}}%
}_{\substack{=a_{s,v_{p}}\\\text{(since }\left(  p+1\right)  -1=p\text{)}%
}}\det\left(  B_{\sim1,\sim\left(  p+1\right)  }\right) \nonumber\\
&  \ \ \ \ \ \ \ \ \ \ \ \ \ \ \ \ \ \ \ \ \left(  \text{here, we have
substituted }p+1\text{ for }q\text{ in the sum}\right) \nonumber\\
&  =a_{s,t}\det\left(  B_{\sim1,\sim1}\right)  +\sum_{p=1}^{k}\left(
-1\right)  ^{p}a_{s,v_{p}}\det\left(  B_{\sim1,\sim\left(  p+1\right)
}\right)  . \label{pf.lem.sol.det.fraser-yeats.1.1}%
\end{align}

Now, we claim that%
\begin{equation}
B_{\sim1,\sim1}=\operatorname*{sub}\nolimits_{u_{1},u_{2},\ldots,u_{k}}%
^{v_{1},v_{2},\ldots,v_{k}}A. \label{pf.lem.sol.det.fraser-yeats.1.B11}%
\end{equation}

[\textit{Proof of (\ref{pf.lem.sol.det.fraser-yeats.1.B11}):} By Definition
\ref{def.submatrix.minor}, we have%
\begin{align*}
B_{\sim1,\sim1}  &  =\operatorname*{sub}\nolimits_{1,2,\ldots,\widehat{1}%
,\ldots,k+1}^{1,2,\ldots,\widehat{1},\ldots,k+1}B\\
&  =\operatorname*{sub}\nolimits_{2,3,\ldots,k+1}^{2,3,\ldots,k+1}%
B\ \ \ \ \ \ \ \ \ \ \left(  \text{since }\left(  1,2,\ldots,\widehat{1}%
,\ldots,k+1\right)  =\left(  2,3,\ldots,k+1\right)  \right) \\
&  =\left(  \underbrace{a_{u_{\left(  x+1\right)  -1},v_{\left(  y+1\right)
-1}}}_{\substack{=a_{u_{x},v_{y}}\\\text{(since }\left(  x+1\right)
-1=x\\\text{and }\left(  y+1\right)  -1=y\text{)}}}\right)  _{1\leq x\leq
k,\ 1\leq y\leq k}\\
&  \ \ \ \ \ \ \ \ \ \ \ \ \ \ \ \ \ \ \ \ \left(  \text{by Definition
\ref{def.submatrix}, since }B=\left(  a_{u_{i-1},v_{j-1}}\right)  _{1\leq
i\leq k+1,\ 1\leq j\leq k+1}\right) \\
&  =\left(  a_{u_{x},v_{y}}\right)  _{1\leq x\leq k,\ 1\leq y\leq k}.
\end{align*}
On the other hand, Definition \ref{def.submatrix} yields%
\[
\operatorname*{sub}\nolimits_{u_{1},u_{2},\ldots,u_{k}}^{v_{1},v_{2}%
,\ldots,v_{k}}A=\left(  a_{u_{x},v_{y}}\right)  _{1\leq x\leq k,\ 1\leq y\leq
k}\ \ \ \ \ \ \ \ \ \ \left(  \text{since }A=\left(  a_{i,j}\right)  _{1\leq
i\leq n,\ 1\leq j\leq n}\right)  .
\]
Comparing these two equalities, we obtain $B_{\sim1,\sim1}=\operatorname*{sub}%
\nolimits_{u_{1},u_{2},\ldots,u_{k}}^{v_{1},v_{2},\ldots,v_{k}}A$. Thus,
(\ref{pf.lem.sol.det.fraser-yeats.1.B11}) is proven.] \medskip

Next, we claim that%
\begin{equation}
B_{\sim1,\sim\left(  p+1\right)  }=\operatorname*{sub}\nolimits_{u_{1}%
,u_{2},\ldots,u_{k}}^{t,v_{1},v_{2},\ldots,\widehat{v_{p}},\ldots,v_{k}}A
\label{pf.lem.sol.det.fraser-yeats.1.B1p1}%
\end{equation}
for each $p\in\left\{  1,2,\ldots,k\right\}  $.

[\textit{Proof of (\ref{pf.lem.sol.det.fraser-yeats.1.B1p1}):} Let
$p\in\left\{  1,2,\ldots,k\right\}  $. Define a map $\phi:\mathbb{Z}%
\rightarrow\mathbb{Z}$ by setting%
\[
\phi\left(  m\right)  =%
\begin{cases}
m, & \text{if }m\leq p;\\
m+1, & \text{if }m>p
\end{cases}
\ \ \ \ \ \ \ \ \ \ \text{for each }m\in\mathbb{Z}.
\]

\begin{vershort}
\noindent Then, we have%
\[
\left(  \phi\left(  1\right)  ,\phi\left(  2\right)  ,\ldots,\phi\left(
k\right)  \right)  =\left(  1,2,\ldots,\widehat{p+1},\ldots,k+1\right)  .
\]

\end{vershort}

\begin{verlong}
\noindent Thus, for each $m\in\left\{  1,2,\ldots,p\right\}  $, we have%
\begin{align}
\phi\left(  m\right)   &  =%
\begin{cases}
m, & \text{if }m\leq p;\\
m+1, & \text{if }m>p
\end{cases}
\nonumber\\
&  =m \label{pf.lem.sol.det.fraser-yeats.1.B1p1.pf.1}%
\end{align}
(since $m\leq p$ (because $m\in\left\{  1,2,\ldots,p\right\}  $)). In other
words,
\[
\left(  \phi\left(  1\right)  ,\phi\left(  2\right)  ,\ldots,\phi\left(
p\right)  \right)  =\left(  1,2,\ldots,p\right)  .
\]
Furthermore, for each $m\in\left\{  p+1,p+2,\ldots,k\right\}  $, we have%
\begin{align}
\phi\left(  m\right)   &  =%
\begin{cases}
m, & \text{if }m\leq p;\\
m+1, & \text{if }m>p
\end{cases}
\ \ \ \ \ \ \ \ \ \ \left(  \text{by the definition of }\phi\right)
\nonumber\\
&  =m+1 \label{pf.lem.sol.det.fraser-yeats.1.B1p1.pf.2}%
\end{align}
(since $m>p$ (because $m\in\left\{  p+1,p+2,\ldots,k\right\}  $)). In other
words,
\begin{align*}
\left(  \phi\left(  p+1\right)  ,\phi\left(  p+2\right)  ,\ldots,\phi\left(
k\right)  \right)   &  =\left(  \left(  p+1\right)  +1,\left(  p+2\right)
+1,\ldots,k+1\right) \\
&  =\left(  p+2,p+3,\ldots,k+1\right)  .
\end{align*}

Now, let us introduce a notation: If $\left(  x_{1},x_{2},\ldots,x_{u}\right)
$ and $\left(  y_{1},y_{2},\ldots,y_{v}\right)  $ are two lists, then we shall
let $\left(  x_{1},x_{2},\ldots,x_{u}\right)  \ast\left(  y_{1},y_{2}%
,\ldots,y_{v}\right)  $ denote the list \newline$\left(  x_{1},x_{2}%
,\ldots,x_{u},y_{1},y_{2},\ldots,y_{v}\right)  $. Then,%
\begin{align*}
&  \left(  \phi\left(  1\right)  ,\phi\left(  2\right)  ,\ldots,\phi\left(
k\right)  \right) \\
&  =\left(  \phi\left(  1\right)  ,\phi\left(  2\right)  ,\ldots,\phi\left(
p\right)  ,\phi\left(  p+1\right)  ,\phi\left(  p+2\right)  ,\ldots
,\phi\left(  k\right)  \right) \\
&  =\underbrace{\left(  \phi\left(  1\right)  ,\phi\left(  2\right)
,\ldots,\phi\left(  p\right)  \right)  }_{=\left(  1,2,\ldots,p\right)  }%
\ast\underbrace{\left(  \phi\left(  p+1\right)  ,\phi\left(  p+2\right)
,\ldots,\phi\left(  k\right)  \right)  }_{=\left(  p+2,p+3,\ldots,k+1\right)
}\\
&  =\left(  1,2,\ldots,p\right)  \ast\left(  p+2,p+3,\ldots,k+1\right) \\
&  =\left(  1,2,\ldots,p,p+2,p+3,\ldots,k+1\right) \\
&  =\left(  1,2,\ldots,\widehat{p+1},\ldots,k+1\right)  .
\end{align*}

\end{verlong}

Now, by Definition \ref{def.submatrix.minor}, we have%
\begin{align}
B_{\sim1,\sim\left(  p+1\right)  }  &  =\operatorname*{sub}%
\nolimits_{1,2,\ldots,\widehat{1},\ldots,k+1}^{1,2,\ldots,\widehat{p+1}%
,\ldots,k+1}B\nonumber\\
&  =\operatorname*{sub}\nolimits_{2,3,\ldots,k+1}^{\phi\left(  1\right)
,\phi\left(  2\right)  ,\ldots,\phi\left(  k\right)  }B\nonumber\\
&  \ \ \ \ \ \ \ \ \ \ \ \ \ \ \ \ \ \ \ \ \left(
\begin{array}
[c]{c}%
\text{because }\left(  1,2,\ldots,\widehat{1},\ldots,k+1\right)  =\left(
2,3,\ldots,k+1\right) \\
\text{and }\left(  1,2,\ldots,\widehat{p+1},\ldots,k+1\right)  =\left(
\phi\left(  1\right)  ,\phi\left(  2\right)  ,\ldots,\phi\left(  k\right)
\right) \\
\text{(since }\left(  \phi\left(  1\right)  ,\phi\left(  2\right)
,\ldots,\phi\left(  k\right)  \right)  =\left(  1,2,\ldots,\widehat{p+1}%
,\ldots,k+1\right)  \text{)}%
\end{array}
\right) \nonumber\\
&  =\left(  \underbrace{a_{u_{\left(  x+1\right)  -1},v_{\phi\left(  y\right)
-1}}}_{\substack{=a_{u_{x},v_{\phi\left(  y\right)  -1}}\\\text{(since
}\left(  x+1\right)  -1=x\text{)}}}\right)  _{1\leq x\leq k,\ 1\leq y\leq
k}\nonumber\\
&  \ \ \ \ \ \ \ \ \ \ \ \ \ \ \ \ \ \ \ \ \left(  \text{by Definition
\ref{def.submatrix}, since }B=\left(  a_{u_{i-1},v_{j-1}}\right)  _{1\leq
i\leq k+1,\ 1\leq j\leq k+1}\right) \nonumber\\
&  =\left(  a_{u_{x},v_{\phi\left(  y\right)  -1}}\right)  _{1\leq x\leq
k,\ 1\leq y\leq k}. \label{pf.lem.sol.det.fraser-yeats.1.B1p1.pf.3}%
\end{align}

\begin{vershort}
On the other hand, because of $t=v_{0}$, we have%
\[
\left(  t,v_{1},v_{2},\ldots,\widehat{v_{p}},\ldots,v_{k}\right)  =\left(
v_{0},v_{1},v_{2},\ldots,\widehat{v_{p}},\ldots,v_{k}\right)  =\left(
v_{\phi\left(  1\right)  -1},v_{\phi\left(  2\right)  -1},\ldots
,v_{\phi\left(  k\right)  -1}\right)
\]
(since the subscripts $0,1,2,\ldots,\widehat{p},\ldots,k$ in the list $\left(
v_{0},v_{1},v_{2},\ldots,\widehat{v_{p}},\ldots,v_{k}\right)  $ are precisely
the numbers $\phi\left(  1\right)  -1,\phi\left(  2\right)  -1,\ldots
,\phi\left(  k\right)  -1$).
\end{vershort}

\begin{verlong}
On the other hand, recall that $t=v_{0}$. Thus,%
\begin{align}
&  \left(  t,v_{1},v_{2},\ldots,\widehat{v_{p}},\ldots,v_{k}\right)
\nonumber\\
&  =\left(  v_{0},v_{1},v_{2},\ldots,\widehat{v_{p}},\ldots,v_{k}\right)
\nonumber\\
&  =\left(  v_{0},v_{1},v_{2},\ldots,v_{p-1},v_{p+1},v_{p+2},\ldots
,v_{k}\right) \nonumber\\
&  =\left(  v_{0},v_{1},v_{2},\ldots,v_{p-1}\right)  \ast\left(
v_{p+1},v_{p+2},\ldots,v_{k}\right)  .
\label{pf.lem.sol.det.fraser-yeats.1.B1p1.pf.4}%
\end{align}

However, for each $m\in\left\{  1,2,\ldots,p\right\}  $, we have $\phi\left(
m\right)  =m$ (by (\ref{pf.lem.sol.det.fraser-yeats.1.B1p1.pf.1})) and
therefore $v_{\phi\left(  m\right)  -1}=v_{m-1}$. Thus,
\[
\left(  v_{\phi\left(  1\right)  -1},v_{\phi\left(  2\right)  -1}%
,\ldots,v_{\phi\left(  p\right)  -1}\right)  =\left(  v_{1-1},v_{2-1}%
,\ldots,v_{p-1}\right)  =\left(  v_{0},v_{1},v_{2},\ldots,v_{p-1}\right)  .
\]

Furthermore, for each $m\in\left\{  p+1,p+2,\ldots,k\right\}  $, we have
$\phi\left(  m\right)  -1=m$ (by
(\ref{pf.lem.sol.det.fraser-yeats.1.B1p1.pf.2})) and therefore $v_{\phi\left(
m\right)  -1}=v_{m}$. Hence,%
\[
\left(  v_{\phi\left(  p+1\right)  -1},v_{\phi\left(  p+2\right)  -1}%
,\ldots,v_{\phi\left(  k\right)  -1}\right)  =\left(  v_{p+1},v_{p+2}%
,\ldots,v_{k}\right)  .
\]
Now,
\begin{align*}
&  \left(  v_{\phi\left(  1\right)  -1},v_{\phi\left(  2\right)  -1}%
,\ldots,v_{\phi\left(  k\right)  -1}\right) \\
&  =\left(  v_{\phi\left(  1\right)  -1},v_{\phi\left(  2\right)  -1}%
,\ldots,v_{\phi\left(  p\right)  -1},v_{\phi\left(  p+1\right)  -1}%
,v_{\phi\left(  p+2\right)  -1},\ldots,v_{\phi\left(  k\right)  -1}\right) \\
&  =\underbrace{\left(  v_{\phi\left(  1\right)  -1},v_{\phi\left(  2\right)
-1},\ldots,v_{\phi\left(  p\right)  -1}\right)  }_{=\left(  v_{0},v_{1}%
,v_{2},\ldots,v_{p-1}\right)  }\ast\underbrace{\left(  v_{\phi\left(
p+1\right)  -1},v_{\phi\left(  p+2\right)  -1},\ldots,v_{\phi\left(  k\right)
-1}\right)  }_{=\left(  v_{p+1},v_{p+2},\ldots,v_{k}\right)  }\\
&  =\left(  v_{0},v_{1},v_{2},\ldots,v_{p-1}\right)  \ast\left(
v_{p+1},v_{p+2},\ldots,v_{k}\right)  .
\end{align*}
Comparing this with (\ref{pf.lem.sol.det.fraser-yeats.1.B1p1.pf.4}), we obtain%
\[
\left(  t,v_{1},v_{2},\ldots,\widehat{v_{p}},\ldots,v_{k}\right)  =\left(
v_{\phi\left(  1\right)  -1},v_{\phi\left(  2\right)  -1},\ldots
,v_{\phi\left(  k\right)  -1}\right)  .
\]

\end{verlong}

Therefore,
\[
\operatorname*{sub}\nolimits_{u_{1},u_{2},\ldots,u_{k}}^{t,v_{1},v_{2}%
,\ldots,\widehat{v_{p}},\ldots,v_{k}}A=\operatorname*{sub}\nolimits_{u_{1}%
,u_{2},\ldots,u_{k}}^{v_{\phi\left(  1\right)  -1},v_{\phi\left(  2\right)
-1},\ldots,v_{\phi\left(  k\right)  -1}}A=\left(  a_{u_{x},v_{\phi\left(
y\right)  -1}}\right)  _{1\leq x\leq k,\ 1\leq y\leq k}%
\]
(by Definition \ref{def.submatrix}, since $A=\left(  a_{i,j}\right)  _{1\leq
i\leq n,\ 1\leq j\leq n}$). Comparing this with
(\ref{pf.lem.sol.det.fraser-yeats.1.B1p1.pf.3}), we find $B_{\sim1,\sim\left(
p+1\right)  }=\operatorname*{sub}\nolimits_{u_{1},u_{2},\ldots,u_{k}}%
^{t,v_{1},v_{2},\ldots,\widehat{v_{p}},\ldots,v_{k}}A$. This proves
(\ref{pf.lem.sol.det.fraser-yeats.1.B1p1}).] \medskip

Now, (\ref{pf.lem.sol.det.fraser-yeats.1.1}) becomes%
\begin{align*}
\det B  &  =a_{s,t}\det\left(  \underbrace{B_{\sim1,\sim1}}%
_{\substack{=\operatorname*{sub}\nolimits_{u_{1},u_{2},\ldots,u_{k}}%
^{v_{1},v_{2},\ldots,v_{k}}A\\\text{(by
(\ref{pf.lem.sol.det.fraser-yeats.1.B11}))}}}\right)  +\sum_{p=1}^{k}\left(
-1\right)  ^{p}a_{s,v_{p}}\det\left(  \underbrace{B_{\sim1,\sim\left(
p+1\right)  }}_{\substack{=\operatorname*{sub}\nolimits_{u_{1},u_{2}%
,\ldots,u_{k}}^{t,v_{1},v_{2},\ldots,\widehat{v_{p}},\ldots,v_{k}}A\\\text{(by
(\ref{pf.lem.sol.det.fraser-yeats.1.B1p1}))}}}\right) \\
&  =a_{s,t}\det\left(  \operatorname*{sub}\nolimits_{u_{1},u_{2},\ldots,u_{k}%
}^{v_{1},v_{2},\ldots,v_{k}}A\right)  +\sum_{p=1}^{k}\left(  -1\right)
^{p}a_{s,v_{p}}\det\left(  \operatorname*{sub}\nolimits_{u_{1},u_{2}%
,\ldots,u_{k}}^{t,v_{1},v_{2},\ldots,\widehat{v_{p}},\ldots,v_{k}}A\right)  .
\end{align*}
In view of $B=\operatorname*{sub}\nolimits_{s,u_{1},u_{2},\ldots,u_{k}%
}^{t,v_{1},v_{2},\ldots,v_{k}}A$, we can rewrite this as%
\begin{align*}
&  \det\left(  \operatorname*{sub}\nolimits_{s,u_{1},u_{2},\ldots,u_{k}%
}^{t,v_{1},v_{2},\ldots,v_{k}}A\right) \\
&  =a_{s,t}\det\left(  \operatorname*{sub}\nolimits_{u_{1},u_{2},\ldots,u_{k}%
}^{v_{1},v_{2},\ldots,v_{k}}A\right)  +\sum_{p=1}^{k}\left(  -1\right)
^{p}a_{s,v_{p}}\det\left(  \operatorname*{sub}\nolimits_{u_{1},u_{2}%
,\ldots,u_{k}}^{t,v_{1},v_{2},\ldots,\widehat{v_{p}},\ldots,v_{k}}A\right)  .
\end{align*}
This proves Lemma \ref{lem.sol.det.fraser-yeats.1}.
\end{proof}

We now step to the solution of Exercise \ref{exe.det.fraser-yeats}.

\begin{proof}
[Solution to Exercise \ref{exe.det.fraser-yeats}.]\textbf{(a)} We shall use
the notation introduced in Definition \ref{def.hat-omit}. Write the $n\times
n$-matrix $A$ in the form $A=\left(  a_{i,j}\right)  _{1\leq i\leq n,\ 1\leq
j\leq n}$. From $t\in\left\{  1,2,\ldots,n\right\}  \setminus V$, we obtain
$t\in\left\{  1,2,\ldots,n\right\}  $ and $t\notin V$.

Let us first observe that%
\begin{equation}
\det\left(  \operatorname*{sub}\nolimits_{s,u_{1},u_{2},\ldots,u_{k}}%
^{t,v_{1},v_{2},\ldots,v_{k}}A\right)  =0 \label{sol.det.fraser-yeats.2}%
\end{equation}
for any $s\in U$.

[\textit{Proof of (\ref{sol.det.fraser-yeats.2}):} Let $s\in U$. Thus, $s\in
U=\left\{  u_{1},u_{2},\ldots,u_{k}\right\}  $. Hence, there exists an
$r\in\left\{  1,2,\ldots,k\right\}  $ such that $s=u_{r}$. Consider this $r$.

\begin{vershort}
From $s=u_{r}$, we see that the $1$-st and the $\left(  r+1\right)  $-st rows
of the matrix $\operatorname*{sub}\nolimits_{s,u_{1},u_{2},\ldots,u_{k}%
}^{t,v_{1},v_{2},\ldots,v_{k}}A$ are equal. Hence, this matrix has two equal
rows. Therefore, Exercise \ref{exe.ps4.6} \textbf{(e)} (applied to $k+1$ and
$\operatorname*{sub}\nolimits_{s,u_{1},u_{2},\ldots,u_{k}}^{t,v_{1}%
,v_{2},\ldots,v_{k}}A$ instead of $n$ and $A$) shows that $\det\left(
\operatorname*{sub}\nolimits_{s,u_{1},u_{2},\ldots,u_{k}}^{t,v_{1}%
,v_{2},\ldots,v_{k}}A\right)  =0$.
\end{vershort}

\begin{verlong}
Note that $r\geq1$, so that $r+1\geq1+1=2>1$. Hence, $r+1\neq1$, so that
$1\neq r+1$.

Set $B=\operatorname*{sub}\nolimits_{s,u_{1},u_{2},\ldots,u_{k}}%
^{t,v_{1},v_{2},\ldots,v_{k}}A$. We extend the $k$-tuple $\left(  u_{1}%
,u_{2},\ldots,u_{k}\right)  $ to a $\left(  k+1\right)  $-tuple $\left(
u_{0},u_{1},\ldots,u_{k}\right)  $ by setting $u_{0}=s$. We extend the
$k$-tuple $\left(  v_{1},v_{2},\ldots,v_{k}\right)  $ to a $\left(
k+1\right)  $-tuple $\left(  v_{0},v_{1},\ldots,v_{k}\right)  $ by setting
$v_{0}=t$. We have%
\[
\left(  u_{0},u_{1},\ldots,u_{k}\right)  =\left(  \underbrace{u_{0}}%
_{=s},u_{1},u_{2},\ldots,u_{k}\right)  =\left(  s,u_{1},u_{2},\ldots
,u_{k}\right)
\]
and%
\[
\left(  v_{0},v_{1},\ldots,v_{k}\right)  =\left(  \underbrace{v_{0}}%
_{=t},v_{1},v_{2},\ldots,v_{k}\right)  =\left(  t,v_{1},v_{2},\ldots
,v_{k}\right)  .
\]
Therefore,%
\[
\operatorname*{sub}\nolimits_{u_{0},u_{1},\ldots,u_{k}}^{v_{0},v_{1}%
,\ldots,v_{k}}A=\operatorname*{sub}\nolimits_{s,u_{1},u_{2},\ldots,u_{k}%
}^{t,v_{1},v_{2},\ldots,v_{k}}A=B.
\]
Hence,
\begin{align*}
B  &  =\operatorname*{sub}\nolimits_{u_{0},u_{1},\ldots,u_{k}}^{v_{0}%
,v_{1},\ldots,v_{k}}A\\
&  =\left(  a_{u_{x-1},v_{y-1}}\right)  _{1\leq x\leq k+1,\ 1\leq y\leq
k+1}\ \ \ \ \ \ \ \ \ \ \left(  \text{by Definition \ref{def.submatrix}, since
}A=\left(  a_{i,j}\right)  _{1\leq i\leq n,\ 1\leq j\leq n}\right)  .
\end{align*}
From this, we obtain both
\begin{align*}
\left(  \text{the }1\text{-st row of the matrix }B\right)   &  =\left(
a_{u_{1-1},v_{y-1}}\right)  _{1\leq x\leq1,\ 1\leq y\leq k+1}\\
&  =\left(  a_{u_{r},v_{y-1}}\right)  _{1\leq x\leq1,\ 1\leq y\leq k+1}%
\end{align*}
(since $u_{1-1}=u_{0}=s=u_{r}$) and%
\begin{align*}
\left(  \text{the }\left(  r+1\right)  \text{-st row of the matrix }B\right)
&  =\left(  a_{u_{\left(  r+1\right)  -1},v_{y-1}}\right)  _{1\leq
x\leq1,\ 1\leq y\leq k+1}\\
&  =\left(  a_{u_{r},v_{y-1}}\right)  _{1\leq x\leq1,\ 1\leq y\leq k+1}%
\end{align*}
(since $\left(  r+1\right)  -1=r$). Comparing these two equalities, we
conclude that%
\[
\left(  \text{the }1\text{-st row of the matrix }B\right)  =\left(  \text{the
}\left(  r+1\right)  \text{-st row of the matrix }B\right)  .
\]
In other words, the $1$-st row and the $\left(  r+1\right)  $-st row of the
matrix $B$ are equal. Thus, the matrix $B$ has two equal rows (since its
$1$-st and its $\left(  r+1\right)  $-st rows are not the same
row\footnote{because $1\neq r+1$}). Therefore, Exercise \ref{exe.ps4.6}
\textbf{(e)} (applied to $k+1$ and $B$ instead of $n$ and $A$) shows that
$\det B=0$. This rewrites as $\det\left(  \operatorname*{sub}%
\nolimits_{s,u_{1},u_{2},\ldots,u_{k}}^{t,v_{1},v_{2},\ldots,v_{k}}A\right)
=0$ (since $B=\operatorname*{sub}\nolimits_{s,u_{1},u_{2},\ldots,u_{k}%
}^{t,v_{1},v_{2},\ldots,v_{k}}A$).
\end{verlong}

Thus, (\ref{sol.det.fraser-yeats.2}) is proved.] \medskip

We define an element $c\in\mathbb{K}$ by%
\begin{equation}
c=\det\left(  \operatorname*{sub}\nolimits_{u_{1},u_{2},\ldots,u_{k}}%
^{v_{1},v_{2},\ldots,v_{k}}A\right)  . \label{sol.det.fraser-yeats.c=}%
\end{equation}
Furthermore, for each $p\in\left\{  1,2,\ldots,k\right\}  $, we define an
element $d_{p}\in\mathbb{K}$ by%
\begin{equation}
d_{p}=\det\left(  \operatorname*{sub}\nolimits_{u_{1},u_{2},\ldots,u_{k}%
}^{t,v_{1},v_{2},\ldots,\widehat{v_{p}},\ldots,v_{k}}A\right)  .
\label{sol.det.fraser-yeats.dp=}%
\end{equation}
Thus, for each $s\in\left\{  1,2,\ldots,n\right\}  $, we have%
\begin{align}
&  \det\left(  \operatorname*{sub}\nolimits_{s,u_{1},u_{2},\ldots,u_{k}%
}^{t,v_{1},v_{2},\ldots,v_{k}}A\right) \nonumber\\
&  =a_{s,t}\underbrace{\det\left(  \operatorname*{sub}\nolimits_{u_{1}%
,u_{2},\ldots,u_{k}}^{v_{1},v_{2},\ldots,v_{k}}A\right)  }%
_{\substack{=c\\\text{(by (\ref{sol.det.fraser-yeats.c=}))}}}+\sum_{p=1}%
^{k}\left(  -1\right)  ^{p}a_{s,v_{p}}\underbrace{\det\left(
\operatorname*{sub}\nolimits_{u_{1},u_{2},\ldots,u_{k}}^{t,v_{1},v_{2}%
,\ldots,\widehat{v_{p}},\ldots,v_{k}}A\right)  }_{\substack{=d_{p}\\\text{(by
(\ref{sol.det.fraser-yeats.dp=}))}}}\nonumber\\
&  \ \ \ \ \ \ \ \ \ \ \ \ \ \ \ \ \ \ \ \ \left(  \text{by Lemma
\ref{lem.sol.det.fraser-yeats.1}}\right) \nonumber\\
&  =a_{s,t}c+\sum_{p=1}^{k}\left(  -1\right)  ^{p}a_{s,v_{p}}d_{p}\nonumber\\
&  =a_{s,t}c+\sum_{z=1}^{k}\left(  -1\right)  ^{z}a_{s,v_{z}}d_{z}
\label{sol.det.fraser-yeats.lem}%
\end{align}
(here, we have renamed the summation index $p$ as $z$).

Now, each $s\in\left\{  1,2,\ldots,n\right\}  $ satisfies either $s\in U$ or
$s\notin U$ (but not both at the same time). Hence,%
\begin{align*}
&  \sum_{s\in\left\{  1,2,\ldots,n\right\}  }\left(  -1\right)  ^{s+t}%
\det\left(  A_{\sim s,\sim t}\right)  \det\left(  \operatorname*{sub}%
\nolimits_{s,u_{1},u_{2},\ldots,u_{k}}^{t,v_{1},v_{2},\ldots,v_{k}}A\right) \\
&  =\sum_{\substack{s\in\left\{  1,2,\ldots,n\right\}  ;\\s\in U}}\left(
-1\right)  ^{s+t}\det\left(  A_{\sim s,\sim t}\right)  \underbrace{\det\left(
\operatorname*{sub}\nolimits_{s,u_{1},u_{2},\ldots,u_{k}}^{t,v_{1}%
,v_{2},\ldots,v_{k}}A\right)  }_{\substack{=0\\\text{(by
(\ref{sol.det.fraser-yeats.2}))}}}\\
&  \ \ \ \ \ \ \ \ \ \ +\underbrace{\sum_{\substack{s\in\left\{
1,2,\ldots,n\right\}  ;\\s\notin U}}}_{=\sum_{s\in\left\{  1,2,\ldots
,n\right\}  \setminus U}}\left(  -1\right)  ^{s+t}\det\left(  A_{\sim s,\sim
t}\right)  \det\left(  \operatorname*{sub}\nolimits_{s,u_{1},u_{2}%
,\ldots,u_{k}}^{t,v_{1},v_{2},\ldots,v_{k}}A\right) \\
&  =\underbrace{\sum_{\substack{s\in\left\{  1,2,\ldots,n\right\}  ;\\s\in
U}}\left(  -1\right)  ^{s+t}\det\left(  A_{\sim s,\sim t}\right)  \cdot0}%
_{=0}\\
&  \ \ \ \ \ \ \ \ \ \ +\sum_{s\in\left\{  1,2,\ldots,n\right\}  \setminus
U}\left(  -1\right)  ^{s+t}\det\left(  A_{\sim s,\sim t}\right)  \det\left(
\operatorname*{sub}\nolimits_{s,u_{1},u_{2},\ldots,u_{k}}^{t,v_{1}%
,v_{2},\ldots,v_{k}}A\right) \\
&  =\sum_{s\in\left\{  1,2,\ldots,n\right\}  \setminus U}\left(  -1\right)
^{s+t}\det\left(  A_{\sim s,\sim t}\right)  \det\left(  \operatorname*{sub}%
\nolimits_{s,u_{1},u_{2},\ldots,u_{k}}^{t,v_{1},v_{2},\ldots,v_{k}}A\right)  .
\end{align*}
Therefore,%
\begin{align*}
&  \sum_{s\in\left\{  1,2,\ldots,n\right\}  \setminus U}\left(  -1\right)
^{s+t}\det\left(  A_{\sim s,\sim t}\right)  \det\left(  \operatorname*{sub}%
\nolimits_{s,u_{1},u_{2},\ldots,u_{k}}^{t,v_{1},v_{2},\ldots,v_{k}}A\right) \\
&  =\underbrace{\sum_{s\in\left\{  1,2,\ldots,n\right\}  }}_{=\sum_{s=1}^{n}%
}\left(  -1\right)  ^{s+t}\underbrace{\det\left(  A_{\sim s,\sim t}\right)
\det\left(  \operatorname*{sub}\nolimits_{s,u_{1},u_{2},\ldots,u_{k}}%
^{t,v_{1},v_{2},\ldots,v_{k}}A\right)  }_{=\det\left(  \operatorname*{sub}%
\nolimits_{s,u_{1},u_{2},\ldots,u_{k}}^{t,v_{1},v_{2},\ldots,v_{k}}A\right)
\det\left(  A_{\sim s,\sim t}\right)  }\\
&  =\sum_{s=1}^{n}\left(  -1\right)  ^{s+t}\underbrace{\det\left(
\operatorname*{sub}\nolimits_{s,u_{1},u_{2},\ldots,u_{k}}^{t,v_{1}%
,v_{2},\ldots,v_{k}}A\right)  }_{\substack{=a_{s,t}c+\sum_{z=1}^{k}\left(
-1\right)  ^{z}a_{s,v_{z}}d_{z}\\\text{(by (\ref{sol.det.fraser-yeats.lem}))}%
}}\det\left(  A_{\sim s,\sim t}\right) \\
&  =\sum_{s=1}^{n}\underbrace{\left(  -1\right)  ^{s+t}\left(  a_{s,t}%
c+\sum_{z=1}^{k}\left(  -1\right)  ^{z}a_{s,v_{z}}d_{z}\right)  \det\left(
A_{\sim s,\sim t}\right)  }_{=\left(  -1\right)  ^{s+t}a_{s,t}c\det\left(
A_{\sim s,\sim t}\right)  +\left(  -1\right)  ^{s+t}\sum_{z=1}^{k}\left(
-1\right)  ^{z}a_{s,v_{z}}d_{z}\det\left(  A_{\sim s,\sim t}\right)  }\\
&  =\sum_{s=1}^{n}\left(  \left(  -1\right)  ^{s+t}a_{s,t}c\det\left(  A_{\sim
s,\sim t}\right)  +\left(  -1\right)  ^{s+t}\sum_{z=1}^{k}\left(  -1\right)
^{z}a_{s,v_{z}}d_{z}\det\left(  A_{\sim s,\sim t}\right)  \right) \\
&  =\sum_{s=1}^{n}\left(  -1\right)  ^{s+t}a_{s,t}c\det\left(  A_{\sim s,\sim
t}\right)  +\sum_{s=1}^{n}\left(  -1\right)  ^{s+t}\sum_{z=1}^{k}\left(
-1\right)  ^{z}a_{s,v_{z}}d_{z}\det\left(  A_{\sim s,\sim t}\right)  .
\end{align*}
In view of%
\begin{align*}
&  \sum_{s=1}^{n}\left(  -1\right)  ^{s+t}a_{s,t}c\det\left(  A_{\sim s,\sim
t}\right) \\
&  =c\sum_{s=1}^{n}\left(  -1\right)  ^{s+t}a_{s,t}\det\left(  A_{\sim s,\sim
t}\right) \\
&  =c\underbrace{\sum_{p=1}^{n}\left(  -1\right)  ^{p+t}a_{p,t}\det\left(
A_{\sim p,\sim t}\right)  }_{\substack{=\det A\\\text{(since Theorem
\ref{thm.laplace.gen} \textbf{(b)} (applied to }q=t\text{)}\\\text{yields
}\det A=\sum_{p=1}^{n}\left(  -1\right)  ^{p+t}a_{p,t}\det\left(  A_{\sim
p,\sim t}\right)  \text{)}}}\\
&  \ \ \ \ \ \ \ \ \ \ \ \ \ \ \ \ \ \ \ \ \left(  \text{here, we have renamed
the summation index }s\text{ as }p\right) \\
&  =c\cdot\det A
\end{align*}
and%
\begin{align*}
&  \sum_{s=1}^{n}\left(  -1\right)  ^{s+t}\sum_{z=1}^{k}\left(  -1\right)
^{z}a_{s,v_{z}}d_{z}\det\left(  A_{\sim s,\sim t}\right) \\
&  =\sum_{z=1}^{k}\left(  -1\right)  ^{z}\underbrace{\sum_{s=1}^{n}\left(
-1\right)  ^{s+t}a_{s,v_{z}}d_{z}\det\left(  A_{\sim s,\sim t}\right)
}_{\substack{=d_{z}\sum_{s=1}^{n}\left(  -1\right)  ^{s+t}a_{s,v_{z}}%
\det\left(  A_{\sim s,\sim t}\right)  \\=d_{z}\sum_{p=1}^{n}\left(  -1\right)
^{p+t}a_{p,v_{z}}\det\left(  A_{\sim p,\sim t}\right)  \\\text{(here, we have
renamed the summation index }s\text{ as }p\text{)}}}\\
&  =\sum_{z=1}^{k}\left(  -1\right)  ^{z}d_{z}\sum_{p=1}^{n}\left(  -1\right)
^{p+t}a_{p,v_{z}}\det\left(  A_{\sim p,\sim t}\right)  ,
\end{align*}
we can rewrite this as%
\begin{align}
&  \sum_{s\in\left\{  1,2,\ldots,n\right\}  \setminus U}\left(  -1\right)
^{s+t}\det\left(  A_{\sim s,\sim t}\right)  \det\left(  \operatorname*{sub}%
\nolimits_{s,u_{1},u_{2},\ldots,u_{k}}^{t,v_{1},v_{2},\ldots,v_{k}}A\right)
\nonumber\\
&  =c\cdot\det A+\sum_{z=1}^{k}\left(  -1\right)  ^{z}d_{z}\sum_{p=1}%
^{n}\left(  -1\right)  ^{p+t}a_{p,v_{z}}\det\left(  A_{\sim p,\sim t}\right)
. \label{sol.det.fraser-yeats.8}%
\end{align}

However, for each $z\in\left\{  1,2,\ldots,k\right\}  $, we have $t\neq v_{z}$
(because otherwise, we would have $t=v_{z}\in\left\{  v_{1},v_{2},\ldots
,v_{k}\right\}  =V$, which would contradict the fact that $t\notin V$) and
therefore%
\begin{equation}
0=\sum_{p=1}^{n}\left(  -1\right)  ^{p+t}a_{p,v_{z}}\det\left(  A_{\sim p,\sim
t}\right)  . \label{sol.det.fraser-yeats.9}%
\end{equation}
(by Proposition \ref{prop.laplace.0} \textbf{(b)}, applied to $t$ and $v_{z}$
instead of $q$ and $r$). Thus, (\ref{sol.det.fraser-yeats.8}) becomes%
\begin{align*}
&  \sum_{s\in\left\{  1,2,\ldots,n\right\}  \setminus U}\left(  -1\right)
^{s+t}\det\left(  A_{\sim s,\sim t}\right)  \det\left(  \operatorname*{sub}%
\nolimits_{s,u_{1},u_{2},\ldots,u_{k}}^{t,v_{1},v_{2},\ldots,v_{k}}A\right) \\
&  =c\cdot\det A+\sum_{z=1}^{k}\left(  -1\right)  ^{z}d_{z}\underbrace{\sum
_{p=1}^{n}\left(  -1\right)  ^{p+t}a_{p,v_{z}}\det\left(  A_{\sim p,\sim
t}\right)  }_{\substack{=0\\\text{(by (\ref{sol.det.fraser-yeats.9}))}}}\\
&  =c\cdot\det A+\underbrace{\sum_{z=1}^{k}\left(  -1\right)  ^{z}d_{z}0}%
_{=0}=c\cdot\det A=\det A\cdot\underbrace{c}_{=\det\left(  \operatorname*{sub}%
\nolimits_{u_{1},u_{2},\ldots,u_{k}}^{v_{1},v_{2},\ldots,v_{k}}A\right)  }\\
&  =\det A\cdot\det\left(  \operatorname*{sub}\nolimits_{u_{1},u_{2}%
,\ldots,u_{k}}^{v_{1},v_{2},\ldots,v_{k}}A\right)  .
\end{align*}
This solves Exercise \ref{exe.det.fraser-yeats} \textbf{(a)}. \medskip

\textbf{(b)} \textit{Alternative proof of Proposition \ref{prop.desnanot.12}.}
Let $n$, $A$, $a_{i,j}$ and $\widetilde{A}$ be as in Proposition
\ref{prop.desnanot.12}. Let $U=\left\{  3,4,\ldots,n\right\}  $ and
$V=\left\{  3,4,\ldots,n\right\}  $. Then, $1\in\left\{  1,2,\ldots,n\right\}
\setminus V$ (since $\left\{  1,2,\ldots,n\right\}  \setminus\underbrace{V}%
_{=\left\{  3,4,\ldots,n\right\}  }=\left\{  1,2,\ldots,n\right\}
\setminus\left\{  3,4,\ldots,n\right\}  =\left\{  1,2\right\}  $ clearly
contains $1$). Hence, Exercise \ref{exe.det.fraser-yeats} \textbf{(a)}
(applied to $k=n-2$ and $\left(  u_{1},u_{2},\ldots,u_{k}\right)  =\left(
3,4,\ldots,n\right)  $ and $\left(  v_{1},v_{2},\ldots,v_{k}\right)  =\left(
3,4,\ldots,n\right)  $ and $t=1$) yields%
\begin{align*}
&  \det A\cdot\det\left(  \operatorname*{sub}\nolimits_{3,4,\ldots
,n}^{3,4,\ldots,n}A\right) \\
&  =\sum_{s\in\left\{  1,2,\ldots,n\right\}  \setminus U}\left(  -1\right)
^{s+1}\det\left(  A_{\sim s,\sim1}\right)  \det\left(  \operatorname*{sub}%
\nolimits_{s,3,4,\ldots,n}^{1,3,4,\ldots,n}A\right) \\
&  =\sum_{s\in\left\{  1,2\right\}  }\left(  -1\right)  ^{s+1}\det\left(
A_{\sim s,\sim1}\right)  \det\left(  \operatorname*{sub}%
\nolimits_{s,3,4,\ldots,n}^{1,3,4,\ldots,n}A\right) \\
&  \ \ \ \ \ \ \ \ \ \ \ \ \ \ \ \ \ \ \ \ \left(  \text{since }\left\{
1,2,\ldots,n\right\}  \setminus\underbrace{U}_{=\left\{  3,4,\ldots,n\right\}
}=\left\{  1,2,\ldots,n\right\}  \setminus\left\{  3,4,\ldots,n\right\}
=\left\{  1,2\right\}  \right) \\
&  =\left(  -1\right)  ^{1+1}\det\left(  A_{\sim1,\sim1}\right)  \det\left(
\operatorname*{sub}\nolimits_{1,3,4,\ldots,n}^{1,3,4,\ldots,n}A\right)
+\left(  -1\right)  ^{2+1}\det\left(  A_{\sim2,\sim1}\right)  \det\left(
\operatorname*{sub}\nolimits_{2,3,4,\ldots,n}^{1,3,4,\ldots,n}A\right)  .
\end{align*}
In view of
\begin{align*}
\operatorname*{sub}\nolimits_{3,4,\ldots,n}^{3,4,\ldots,n}A  &  =\left(
a_{x+2,y+2}\right)  _{1\leq x\leq n-2,\ 1\leq y\leq n-2}%
\ \ \ \ \ \ \ \ \ \ \left(
\begin{array}
[c]{c}%
\text{by Definition \ref{def.submatrix},}\\
\text{since }A=\left(  a_{i,j}\right)  _{1\leq i\leq n,\ 1\leq j\leq n}%
\end{array}
\right) \\
&  =\left(  a_{i+2,j+2}\right)  _{1\leq i\leq n-2,\ 1\leq j\leq n-2}%
\ \ \ \ \ \ \ \ \ \ \left(
\begin{array}
[c]{c}%
\text{here, we have renamed the}\\
\text{index }\left(  x,y\right)  \text{ as }\left(  i,j\right)
\end{array}
\right) \\
&  =\widetilde{A}\ \ \ \ \ \ \ \ \ \ \left(  \text{since }\widetilde{A}\text{
is defined as }\left(  a_{i+2,j+2}\right)  _{1\leq i\leq n-2,\ 1\leq j\leq
n-2}\right)
\end{align*}
and%
\begin{align*}
\operatorname*{sub}\nolimits_{1,3,4,\ldots,n}^{1,3,4,\ldots,n}A  &
=\operatorname*{sub}\nolimits_{1,2,\ldots,\widehat{2},\ldots,n}^{1,2,\ldots
,\widehat{2},\ldots,n}A\ \ \ \ \ \ \ \ \ \ \left(  \text{since }\left(
1,3,4,\ldots,n\right)  =\left(  1,2,\ldots,\widehat{2},\ldots,n\right)
\right) \\
&  =A_{\sim2,\sim2}\ \ \ \ \ \ \ \ \ \ \left(  \text{since }A_{\sim2,\sim
2}\text{ is defined as }\operatorname*{sub}\nolimits_{1,2,\ldots
,\widehat{2},\ldots,n}^{1,2,\ldots,\widehat{2},\ldots,n}A\right)
\end{align*}
and%
\begin{align*}
\operatorname*{sub}\nolimits_{2,3,4,\ldots,n}^{1,3,4,\ldots,n}A  &
=\operatorname*{sub}\nolimits_{1,2,\ldots,\widehat{1},\ldots,n}^{1,2,\ldots
,\widehat{2},\ldots,n}A\ \ \ \ \ \ \ \ \ \ \left(
\begin{array}
[c]{c}%
\text{since }\left(  1,3,4,\ldots,n\right)  =\left(  1,2,\ldots,\widehat{2}%
,\ldots,n\right) \\
\text{and }\left(  2,3,4,\ldots,n\right)  =\left(  1,2,\ldots,\widehat{1}%
,\ldots,n\right)
\end{array}
\right) \\
&  =A_{\sim1,\sim2}\ \ \ \ \ \ \ \ \ \ \left(  \text{since }A_{\sim1,\sim
2}\text{ is defined as }\operatorname*{sub}\nolimits_{1,2,\ldots
,\widehat{1},\ldots,n}^{1,2,\ldots,\widehat{2},\ldots,n}A\right)  ,
\end{align*}
we can rewrite this equality as%
\begin{align*}
&  \det A\cdot\det\widetilde{A}\\
&  =\underbrace{\left(  -1\right)  ^{1+1}}_{=\left(  -1\right)  ^{2}%
=1}\underbrace{\det\left(  A_{\sim1,\sim1}\right)  \det\left(  A_{\sim2,\sim
2}\right)  }_{=\det\left(  A_{\sim1,\sim1}\right)  \cdot\det\left(
A_{\sim2,\sim2}\right)  }+\underbrace{\left(  -1\right)  ^{2+1}}_{=\left(
-1\right)  ^{3}=-1}\underbrace{\det\left(  A_{\sim2,\sim1}\right)  \det\left(
A_{\sim1,\sim2}\right)  }_{=\det\left(  A_{\sim1,\sim2}\right)  \cdot
\det\left(  A_{\sim2,\sim1}\right)  }\\
&  =\det\left(  A_{\sim1,\sim1}\right)  \cdot\det\left(  A_{\sim2,\sim
2}\right)  -\det\left(  A_{\sim1,\sim2}\right)  \cdot\det\left(  A_{\sim
2,\sim1}\right)  .
\end{align*}
Thus, Proposition \ref{prop.desnanot.12} is proved once again.

[\textit{Remark:} We could now easily rederive Proposition
\ref{prop.desnanot.1n} and the more general Theorem \ref{thm.desnanot} from
Proposition \ref{prop.desnanot.12} by appropriately permuting rows and
columns. However, this boring exercise is best left to the reader.]
\end{proof}

\subsection{Solution to Exercise \ref{exe.det.syl-lem}}

In this section, we shall use the notations introduced in Definition
\ref{def.sect.laplace.notations}. We begin with two lemmas:

\begin{lemma}
\label{lem.sol.det.syl-lem.0}Let $n\in\mathbb{N}$ and $m\in\mathbb{N}$ and
$p\in\mathbb{N}$. Let $A\in\mathbb{K}^{n\times m}$ and $B\in\mathbb{K}%
^{n\times p}$ be two matrices. Let $J$ be a subset of $\left\{  1,2,\ldots
,m\right\}  $, and let $K$ be a subset of $\left\{  1,2,\ldots,p\right\}  $
such that $\left\vert J\right\vert =\left\vert K\right\vert $. Define an
$n\times m$-matrix $C\in\mathbb{K}^{n\times m}$ by%
\begin{equation}
C=\left(
\begin{array}
[c]{c}%
A\leftarrow B\\
J\leftarrow K
\end{array}
\right)  . \label{eq.lem.sol.det.syl-lem.0.C=}%
\end{equation}
Let $R$ be a subset of $\left\{  1,2,\ldots,n\right\}  $. Then:

\textbf{(a)} We have $\operatorname*{sub}\nolimits_{w\left(  R\right)
}^{w\left(  J\right)  }C=\operatorname*{sub}\nolimits_{w\left(  R\right)
}^{w\left(  K\right)  }B$.

\textbf{(b)} We have $\operatorname*{sub}\nolimits_{w\left(  R\right)
}^{w\left(  \widetilde{J}\right)  }C=\operatorname*{sub}\nolimits_{w\left(
R\right)  }^{w\left(  \widetilde{J}\right)  }A$, where $\widetilde{J}$ denotes
the complement $\left\{  1,2,\ldots,m\right\}  \setminus J$ of $J$.
\end{lemma}

\begin{proof}
[Proof of Lemma \ref{lem.sol.det.syl-lem.0}.]Roughly speaking, part
\textbf{(a)} of this lemma is saying that a bunch of columns of $C$ (namely,
its $j$-th columns for $j\in J$) are identical with certain columns of $B$,
whereas part \textbf{(b)} is saying that the remaining columns of $C$ are
identical with the respective columns of $A$. Both of these claims are obvious
consequences of how the matrix $C$ was defined. For the sake of completeness,
let us nevertheless give a rigorous proof that relies on the formal definition
of $C$.

Write the matrices $A$ and $B$ as $A=\left(  a_{x,y}\right)  _{1\leq x\leq
n,\ 1\leq y\leq m}$ and $B=\left(  b_{x,y}\right)  _{1\leq x\leq n,\ 1\leq
y\leq p}$.

Let $j_{1},j_{2},\ldots,j_{q}$ be all elements of $J$, listed in increasing
order (with no repetitions). Thus, $q=\left\vert J\right\vert =\left\vert
K\right\vert $. Hence, the set $K$ has exactly $q$ elements. Let $k_{1}%
,k_{2},\ldots,k_{q}$ be all elements of $K$, listed in increasing order (with
no repetitions). (This is well-defined, since $K$ has exactly $q$ elements.)

For any $x\in\left\{  1,2,\ldots,n\right\}  $ and $y\in\left\{  1,2,\ldots
,m\right\}  $, we define an element $c_{x,y}\in\mathbb{K}$ by the equality
(\ref{eq.def.matrix-col-exchange.cxy=}). Then, $\left(
\begin{array}
[c]{c}%
A\leftarrow B\\
J\leftarrow K
\end{array}
\right)  $ is the $n\times m$-matrix $\left(  c_{x,y}\right)  _{1\leq x\leq
n,\ 1\leq y\leq m}$ (because of how we defined $\left(
\begin{array}
[c]{c}%
A\leftarrow B\\
J\leftarrow K
\end{array}
\right)  $ in Definition \ref{def.matrix-col-exchange}). Thus, we have
\[
\left(
\begin{array}
[c]{c}%
A\leftarrow B\\
J\leftarrow K
\end{array}
\right)  =\left(  c_{x,y}\right)  _{1\leq x\leq n,\ 1\leq y\leq m}.
\]
In view of (\ref{eq.lem.sol.det.syl-lem.0.C=}), this rewrites as
\begin{equation}
C=\left(  c_{x,y}\right)  _{1\leq x\leq n,\ 1\leq y\leq m}.
\label{pf.lem.sol.det.syl-lem.0.C=}%
\end{equation}

Recall that $j_{1},j_{2},\ldots,j_{q}$ are all elements of $J$, listed in
increasing order (with no repetitions). In other words, $\left(  j_{1}%
,j_{2},\ldots,j_{q}\right)  $ is the list of all elements of $J$ in increasing
order (with no repetitions). In other words, $\left(  j_{1},j_{2},\ldots
,j_{q}\right)  =w\left(  J\right)  $ (since $w\left(  J\right)  $ is defined
to be the list of all elements of $J$ in increasing order (with no
repetitions)). Likewise, we see that $\left(  k_{1},k_{2},\ldots,k_{q}\right)
=w\left(  K\right)  $. Let $\left(  r_{1},r_{2},\ldots,r_{s}\right)  $ be the
list $w\left(  R\right)  $. Let $\left(  \ell_{1},\ell_{2},\ldots,\ell
_{t}\right)  $ be the list $w\left(  \widetilde{J}\right)  $. \medskip

\textbf{(a)} We first claim the following:

\begin{statement}
\textit{Claim 1:} We have $c_{r_{x},j_{y}}=b_{r_{x},k_{y}}$ for each
$x\in\left\{  1,2,\ldots,s\right\}  $ and $y\in\left\{  1,2,\ldots,q\right\}
$.
\end{statement}

[\textit{Proof of Claim 1:} Let $x\in\left\{  1,2,\ldots,s\right\}  $ and
$y\in\left\{  1,2,\ldots,q\right\}  $. Then, $j_{y}=j_{i}$ for some
$i\in\left\{  1,2,\ldots,q\right\}  $ (namely, for $i=y$). Moreover, we have
$j_{y}\in\left\{  j_{1},j_{2},\ldots,j_{q}\right\}  =J$ (since $j_{1}%
,j_{2},\ldots,j_{q}$ are all elements of $J$); thus, we don't have
$j_{y}\notin J$. However, (\ref{eq.def.matrix-col-exchange.cxy=}) (applied to
$r_{x}$ and $j_{y}$ instead of $x$ and $y$) yields%
\[
c_{r_{x},j_{y}}=%
\begin{cases}
a_{r_{x},j_{y}}, & \text{if }j_{y}\notin J;\\
b_{r_{x},k_{i}}, & \text{if }j_{y}=j_{i}\text{ for some }i\in\left\{
1,2,\ldots,q\right\}
\end{cases}
.
\]
Since we don't have $j_{y}\notin J$, we can simplify this as follows:%
\[
c_{r,j_{y}}=b_{r_{x},k_{i}},\ \ \ \ \ \ \ \ \ \ \text{if }j_{y}=j_{i}\text{
for some }i\in\left\{  1,2,\ldots,q\right\}  .
\]
We can apply this to $i=y$ (because we have $j_{y}=j_{i}$ for $i=y$). Thus, we
obtain $c_{x,j_{y}}=b_{r_{x},k_{y}}$. This proves Claim 1.] \medskip

Now, recall that $w\left(  R\right)  =\left(  r_{1},r_{2},\ldots,r_{s}\right)
$ (by the definition of $\left(  r_{1},r_{2},\ldots,r_{s}\right)  $) and
$w\left(  J\right)  =\left(  j_{1},j_{2},\ldots,j_{q}\right)  $. Hence,%
\begin{align}
\operatorname*{sub}\nolimits_{w\left(  R\right)  }^{w\left(  J\right)  }C  &
=\operatorname*{sub}\nolimits_{\left(  r_{1},r_{2},\ldots,r_{s}\right)
}^{\left(  j_{1},j_{2},\ldots,j_{q}\right)  }C=\operatorname*{sub}%
\nolimits_{r_{1},r_{2},\ldots,r_{s}}^{j_{1},j_{2},\ldots,j_{q}}C\nonumber\\
&  =\left(  \underbrace{c_{r_{x},j_{y}}}_{\substack{=b_{r_{x},k_{y}%
}\\\text{(by Claim 1)}}}\right)  _{1\leq x\leq s,\ 1\leq y\leq q}%
\ \ \ \ \ \ \ \ \ \ \left(
\begin{array}
[c]{c}%
\text{by Definition \ref{def.submatrix},}\\
\text{since }C=\left(  c_{x,y}\right)  _{1\leq x\leq n,\ 1\leq y\leq m}%
\end{array}
\right) \nonumber\\
&  =\left(  b_{r_{x},k_{y}}\right)  _{1\leq x\leq s,\ 1\leq y\leq q}.
\label{pf.lem.sol.det.syl-lem.0.a.6}%
\end{align}
On the other hand,%
\begin{align*}
\operatorname*{sub}\nolimits_{w\left(  R\right)  }^{w\left(  K\right)  }B  &
=\operatorname*{sub}\nolimits_{\left(  r_{1},r_{2},\ldots,r_{s}\right)
}^{\left(  k_{1},k_{2},\ldots,k_{q}\right)  }B\ \ \ \ \ \ \ \ \ \ \left(
\begin{array}
[c]{c}%
\text{since }w\left(  K\right)  =\left(  k_{1},k_{2},\ldots,k_{q}\right) \\
\text{and }w\left(  R\right)  =\left(  r_{1},r_{2},\ldots,r_{s}\right)
\end{array}
\right) \\
&  =\operatorname*{sub}\nolimits_{r_{1},r_{2},\ldots,r_{s}}^{k_{1}%
,k_{2},\ldots,k_{q}}B=\left(  b_{r_{x},k_{y}}\right)  _{1\leq x\leq s,\ 1\leq
y\leq q}%
\end{align*}
(by Definition \ref{def.submatrix}, since $B=\left(  b_{x,y}\right)  _{1\leq
x\leq n,\ 1\leq y\leq p}$). Comparing this with
(\ref{pf.lem.sol.det.syl-lem.0.a.6}), we obtain $\operatorname*{sub}%
\nolimits_{w\left(  R\right)  }^{w\left(  J\right)  }C=\operatorname*{sub}%
\nolimits_{w\left(  R\right)  }^{w\left(  K\right)  }B$. This proves Lemma
\ref{lem.sol.det.syl-lem.0} \textbf{(a)}. \medskip

\textbf{(b)} Let $\widetilde{J}$ denote the complement $\left\{
1,2,\ldots,m\right\}  \setminus J$ of $J$. We first claim the following:

\begin{statement}
\textit{Claim 2:} We have $c_{r_{x},\ell_{y}}=a_{r_{x},\ell_{y}}$ for each
$x\in\left\{  1,2,\ldots,s\right\}  $ and $y\in\left\{  1,2,\ldots,t\right\}
$.
\end{statement}

[\textit{Proof of Claim 2:} Let $x\in\left\{  1,2,\ldots,s\right\}  $ and
$y\in\left\{  1,2,\ldots,q\right\}  $. Note that $w\left(  \widetilde{J}%
\right)  $ is the list of all elements of $\widetilde{J}$ in increasing order
(with no repetitions) (by the definition of $w\left(  \widetilde{J}\right)
$). Hence, each entry of $w\left(  \widetilde{J}\right)  $ is an element of
$\widetilde{J}$. In other words, each entry of $\left(  \ell_{1},\ell
_{2},\ldots,\ell_{t}\right)  $ is an element of $\widetilde{J}$ (since
$\left(  \ell_{1},\ell_{2},\ldots,\ell_{t}\right)  =w\left(  \widetilde{J}%
\right)  $). Therefore, $\ell_{y}$ is an element of $\widetilde{J}$ (since
$\ell_{y}$ is an entry of $\left(  \ell_{1},\ell_{2},\ldots,\ell_{t}\right)
$). In other words, $\ell_{y}\in\widetilde{J}=\left\{  1,2,\ldots,m\right\}
\setminus J$. Therefore, $\ell_{y}\notin J$.

Now, (\ref{eq.def.matrix-col-exchange.cxy=}) (applied to $r_{x}$ and $\ell
_{y}$ instead of $x$ and $y$) yields%
\begin{align*}
c_{r_{x},\ell_{y}}  &  =%
\begin{cases}
a_{r_{x},\ell_{y}}, & \text{if }\ell_{y}\notin J;\\
b_{r_{x},k_{i}}, & \text{if }\ell_{y}=j_{i}\text{ for some }i\in\left\{
1,2,\ldots,q\right\}
\end{cases}
\\
&  =a_{r_{x},\ell_{y}}\ \ \ \ \ \ \ \ \ \ \left(  \text{since }\ell_{y}\notin
J\right)  .
\end{align*}
This proves Claim 2.] \medskip

Now, recall that $w\left(  R\right)  =\left(  r_{1},r_{2},\ldots,r_{s}\right)
$ (by the definition of $\left(  r_{1},r_{2},\ldots,r_{s}\right)  $) and
$w\left(  \widetilde{J}\right)  =\left(  \ell_{1},\ell_{2},\ldots,\ell
_{t}\right)  $ (by the definition of $\left(  \ell_{1},\ell_{2},\ldots
,\ell_{t}\right)  $). Hence,%
\begin{align}
\operatorname*{sub}\nolimits_{w\left(  R\right)  }^{w\left(  \widetilde{J}%
\right)  }C  &  =\operatorname*{sub}\nolimits_{\left(  r_{1},r_{2}%
,\ldots,r_{s}\right)  }^{\left(  \ell_{1},\ell_{2},\ldots,\ell_{t}\right)
}C=\operatorname*{sub}\nolimits_{r_{1},r_{2},\ldots,r_{s}}^{\ell_{1},\ell
_{2},\ldots,\ell_{t}}C\nonumber\\
&  =\left(  \underbrace{c_{r_{x},\ell_{y}}}_{\substack{=a_{r_{x},\ell_{y}%
}\\\text{(by Claim 2)}}}\right)  _{1\leq x\leq s,\ 1\leq y\leq t}%
\ \ \ \ \ \ \ \ \ \ \left(
\begin{array}
[c]{c}%
\text{by Definition \ref{def.submatrix},}\\
\text{since }C=\left(  c_{x,y}\right)  _{1\leq x\leq n,\ 1\leq y\leq m}%
\end{array}
\right) \nonumber\\
&  =\left(  a_{r_{x},\ell_{y}}\right)  _{1\leq x\leq s,\ 1\leq y\leq t}.
\label{pf.lem.sol.det.syl-lem.0.b.6}%
\end{align}
On the other hand,%
\begin{align*}
\operatorname*{sub}\nolimits_{w\left(  R\right)  }^{w\left(  \widetilde{J}%
\right)  }A  &  =\operatorname*{sub}\nolimits_{\left(  r_{1},r_{2}%
,\ldots,r_{s}\right)  }^{\left(  \ell_{1},\ell_{2},\ldots,\ell_{t}\right)
}A\ \ \ \ \ \ \ \ \ \ \left(
\begin{array}
[c]{c}%
\text{since }w\left(  \widetilde{J}\right)  =\left(  \ell_{1},\ell_{2}%
,\ldots,\ell_{t}\right) \\
\text{and }w\left(  R\right)  =\left(  r_{1},r_{2},\ldots,r_{s}\right)
\end{array}
\right) \\
&  =\operatorname*{sub}\nolimits_{r_{1},r_{2},\ldots,r_{s}}^{\ell_{1},\ell
_{2},\ldots,\ell_{t}}A=\left(  a_{r_{x},\ell_{y}}\right)  _{1\leq x\leq
s,\ 1\leq y\leq t}%
\end{align*}
(by Definition \ref{def.submatrix}, since $A=\left(  a_{x,y}\right)  _{1\leq
x\leq n,\ 1\leq y\leq m}$). Comparing this with
(\ref{pf.lem.sol.det.syl-lem.0.b.6}), we obtain $\operatorname*{sub}%
\nolimits_{w\left(  R\right)  }^{w\left(  \widetilde{J}\right)  }%
C=\operatorname*{sub}\nolimits_{w\left(  R\right)  }^{w\left(  \widetilde{J}%
\right)  }A$. This proves Lemma \ref{lem.sol.det.syl-lem.0} \textbf{(b)}.
\end{proof}

\begin{lemma}
\label{lem.sol.det.syl-lem.1}Let $n\in\mathbb{N}$. Let $A\in\mathbb{K}%
^{n\times n}$ and $B\in\mathbb{K}^{n\times n}$ be two $n\times n$-matrices.
For any subset $I$ of $\left\{  1,2,\ldots,n\right\}  $, we let $\widetilde{I}%
$ denote the complement $\left\{  1,2,\ldots,n\right\}  \setminus I$ of $I$.
(For instance, if $n=4$ and $I=\left\{  1,4\right\}  $, then $\widetilde{I}%
=\left\{  2,3\right\}  $.)

Let $J$ and $K$ be two subsets of $\left\{  1,2,\ldots,n\right\}  $ such that
$\left\vert K\right\vert =\left\vert J\right\vert $. Then,%
\[
\det\left(
\begin{array}
[c]{c}%
A\leftarrow B\\
J\leftarrow K
\end{array}
\right)  =\sum_{\substack{P\subseteq\left\{  1,2,\ldots,n\right\}
;\\\left\vert P\right\vert =\left\vert J\right\vert }}\left(  -1\right)
^{\sum P+\sum J}\det\left(  \operatorname*{sub}\nolimits_{w\left(  P\right)
}^{w\left(  K\right)  }B\right)  \det\left(  \operatorname*{sub}%
\nolimits_{w\left(  \widetilde{P}\right)  }^{w\left(  \widetilde{J}\right)
}A\right)  .
\]

\end{lemma}

\begin{proof}
[Proof of Lemma \ref{lem.sol.det.syl-lem.1}.]Define an $n\times n$-matrix
$C\in\mathbb{K}^{n\times n}$ by%
\begin{equation}
C=\left(
\begin{array}
[c]{c}%
A\leftarrow B\\
J\leftarrow K
\end{array}
\right)  . \label{pf.lem.sol.det.syl-lem.1.C=}%
\end{equation}
Then, Theorem \ref{thm.det.laplace-multi} \textbf{(b)} (applied to $C$ and $J$
instead of $A$ and $Q$) yields%
\[
\det C=\sum_{\substack{P\subseteq\left\{  1,2,\ldots,n\right\}  ;\\\left\vert
P\right\vert =\left\vert J\right\vert }}\left(  -1\right)  ^{\sum P+\sum
J}\det\left(  \operatorname*{sub}\nolimits_{w\left(  P\right)  }^{w\left(
J\right)  }C\right)  \det\left(  \operatorname*{sub}\nolimits_{w\left(
\widetilde{P}\right)  }^{w\left(  \widetilde{J}\right)  }C\right)  .
\]
However, for each subset $P$ of $\left\{  1,2,\ldots,n\right\}  $, we have%
\begin{equation}
\operatorname*{sub}\nolimits_{w\left(  P\right)  }^{w\left(  J\right)
}C=\operatorname*{sub}\nolimits_{w\left(  P\right)  }^{w\left(  K\right)  }B
\label{pf.lem.sol.det.syl-lem.1.2}%
\end{equation}
(by Lemma \ref{lem.sol.det.syl-lem.0} \textbf{(a)}, applied to $m=n$ and $p=n$
and $R=P$) and%
\begin{equation}
\operatorname*{sub}\nolimits_{w\left(  \widetilde{P}\right)  }^{w\left(
\widetilde{J}\right)  }C=\operatorname*{sub}\nolimits_{w\left(  \widetilde{P}%
\right)  }^{w\left(  \widetilde{J}\right)  }A
\label{pf.lem.sol.det.syl-lem.1.3}%
\end{equation}
(by Lemma \ref{lem.sol.det.syl-lem.0} \textbf{(b)}, applied to $m=n$ and $p=n$
and $R=\widetilde{P}$), since $\widetilde{J}$ is the complement $\left\{
1,2,\ldots,n\right\}  \setminus J$ of $J$. Thus,
\begin{align*}
\det C  &  =\sum_{\substack{P\subseteq\left\{  1,2,\ldots,n\right\}
;\\\left\vert P\right\vert =\left\vert J\right\vert }}\left(  -1\right)
^{\sum P+\sum J}\det\underbrace{\left(  \operatorname*{sub}\nolimits_{w\left(
P\right)  }^{w\left(  J\right)  }C\right)  }_{\substack{=\operatorname*{sub}%
\nolimits_{w\left(  P\right)  }^{w\left(  K\right)  }B\\\text{(by
(\ref{pf.lem.sol.det.syl-lem.1.2}))}}}\det\underbrace{\left(
\operatorname*{sub}\nolimits_{w\left(  \widetilde{P}\right)  }^{w\left(
\widetilde{J}\right)  }C\right)  }_{\substack{=\operatorname*{sub}%
\nolimits_{w\left(  \widetilde{P}\right)  }^{w\left(  \widetilde{J}\right)
}A\\\text{(by (\ref{pf.lem.sol.det.syl-lem.1.3}))}}}\\
&  =\sum_{\substack{P\subseteq\left\{  1,2,\ldots,n\right\}  ;\\\left\vert
P\right\vert =\left\vert J\right\vert }}\left(  -1\right)  ^{\sum P+\sum
J}\det\left(  \operatorname*{sub}\nolimits_{w\left(  P\right)  }^{w\left(
K\right)  }B\right)  \det\left(  \operatorname*{sub}\nolimits_{w\left(
\widetilde{P}\right)  }^{w\left(  \widetilde{J}\right)  }A\right)  .
\end{align*}
In view of (\ref{pf.lem.sol.det.syl-lem.1.C=}), we can rewrite this as
\[
\det\left(
\begin{array}
[c]{c}%
A\leftarrow B\\
J\leftarrow K
\end{array}
\right)  =\sum_{\substack{P\subseteq\left\{  1,2,\ldots,n\right\}
;\\\left\vert P\right\vert =\left\vert J\right\vert }}\left(  -1\right)
^{\sum P+\sum J}\det\left(  \operatorname*{sub}\nolimits_{w\left(  P\right)
}^{w\left(  K\right)  }B\right)  \det\left(  \operatorname*{sub}%
\nolimits_{w\left(  \widetilde{P}\right)  }^{w\left(  \widetilde{J}\right)
}A\right)  .
\]
This proves Lemma \ref{lem.sol.det.syl-lem.1}.
\end{proof}

We are now ready to solve Exercise \ref{exe.det.syl-lem}:

\begin{proof}
[Solution to Exercise \ref{exe.det.syl-lem}.]For any subset $I$ of $\left\{
1,2,\ldots,n\right\}  $, we let $\widetilde{I}$ denote the complement
$\left\{  1,2,\ldots,n\right\}  \setminus I$ of $I$. (For instance, if $n=4$
and $I=\left\{  1,4\right\}  $, then $\widetilde{I}=\left\{  2,3\right\}  $.)

If $P$ and $R$ are two subsets of $\left\{  1,2,\ldots,n\right\}  $ satisfying
$\left\vert P\right\vert =\left\vert J\right\vert $ and $\left\vert
R\right\vert =\left\vert J\right\vert $, then we define an element
$\alpha_{P,R}\in\mathbb{K}$ by%
\begin{equation}
\alpha_{P,R}=\left(  -1\right)  ^{\sum R+\sum J}\det\left(
\operatorname*{sub}\nolimits_{w\left(  P\right)  }^{w\left(  J\right)
}A\right)  \det\left(  \operatorname*{sub}\nolimits_{w\left(  \widetilde{R}%
\right)  }^{w\left(  \widetilde{J}\right)  }A\right)  .
\label{sol.det.syl-lem.beta=}%
\end{equation}
(This is easily seen to be well-defined, because both $\operatorname*{sub}%
\nolimits_{w\left(  P\right)  }^{w\left(  J\right)  }A$ and
$\operatorname*{sub}\nolimits_{w\left(  \widetilde{R}\right)  }^{w\left(
\widetilde{J}\right)  }A$ are square matrices\footnote{\textit{Proof.} Let $P$
and $R$ be two subsets of $\left\{  1,2,\ldots,n\right\}  $ satisfying
$\left\vert P\right\vert =\left\vert J\right\vert $ and $\left\vert
R\right\vert =\left\vert J\right\vert $. The definition of $\widetilde{R}$
yields $\widetilde{R}=\left\{  1,2,\ldots,n\right\}  \setminus R$, so that
\begin{align*}
\left\vert \widetilde{R}\right\vert  &  =\left\vert \left\{  1,2,\ldots
,n\right\}  \setminus R\right\vert =\underbrace{\left\vert \left\{
1,2,\ldots,n\right\}  \right\vert }_{=n}-\left\vert R\right\vert
\ \ \ \ \ \ \ \ \ \ \left(  \text{since }R\subseteq\left\{  1,2,\ldots
,n\right\}  \right) \\
&  =n-\left\vert R\right\vert .
\end{align*}
Similarly, $\left\vert \widetilde{J}\right\vert =n-\left\vert J\right\vert $.
Now, $\left\vert \widetilde{R}\right\vert =n-\underbrace{\left\vert
R\right\vert }_{=\left\vert J\right\vert }=n-\left\vert J\right\vert
=\left\vert \widetilde{J}\right\vert $.
\par
Now, $w\left(  J\right)  $ is a $\left\vert J\right\vert $-tuple of elements
of $J$ (by Proposition \ref{prop.sect.laplace.notations.w(I)}, applied to
$I=J$). Thus, in particular, $w\left(  J\right)  $ is a $\left\vert
J\right\vert $-tuple. Similarly, $w\left(  P\right)  $ is a $\left\vert
P\right\vert $-tuple. In other words, $w\left(  P\right)  $ is a $\left\vert
J\right\vert $-tuple (since $\left\vert P\right\vert =\left\vert J\right\vert
$). Now, $\operatorname*{sub}\nolimits_{w\left(  P\right)  }^{w\left(
J\right)  }A$ is a $\left\vert J\right\vert \times\left\vert J\right\vert
$-matrix (since $w\left(  J\right)  $ and $w\left(  P\right)  $ are
$\left\vert J\right\vert $-tuples), hence a square matrix.
\par
It remains to show that $\operatorname*{sub}\nolimits_{w\left(  \widetilde{R}%
\right)  }^{w\left(  \widetilde{J}\right)  }A$ is a square matrix. However,
$w\left(  \widetilde{J}\right)  $ is a $\widetilde{J}$-tuple of elements of
$\widetilde{J}$ (by Proposition \ref{prop.sect.laplace.notations.w(I)},
applied to $I=\widetilde{J}$). Thus, in particular, $w\left(  \widetilde{J}%
\right)  $ is a $\left\vert \widetilde{J}\right\vert $-tuple. Similarly,
$w\left(  \widetilde{R}\right)  $ is a $\left\vert \widetilde{R}\right\vert
$-tuple. In other words, $w\left(  \widetilde{R}\right)  $ is a $\left\vert
\widetilde{J}\right\vert $-tuple (since $\left\vert \widetilde{R}\right\vert
=\left\vert \widetilde{J}\right\vert $). Now, $\operatorname*{sub}%
\nolimits_{w\left(  \widetilde{R}\right)  }^{w\left(  \widetilde{J}\right)
}A$ is a $\left\vert \widetilde{J}\right\vert \times\left\vert \widetilde{J}%
\right\vert $-matrix (since $w\left(  \widetilde{J}\right)  $ and $w\left(
\widetilde{R}\right)  $ are $\left\vert \widetilde{J}\right\vert $-tuples),
hence a square matrix.
\par
Thus, we have shown that both $\operatorname*{sub}\nolimits_{w\left(
P\right)  }^{w\left(  J\right)  }A$ and $\operatorname*{sub}%
\nolimits_{w\left(  \widetilde{R}\right)  }^{w\left(  \widetilde{J}\right)
}A$ are square matrices.}.)

Let $K$ be a subset of $\left\{  1,2,\ldots,n\right\}  $ satisfying
$\left\vert K\right\vert =\left\vert J\right\vert $. Then, Lemma
\ref{lem.sol.det.syl-lem.1} yields%
\begin{align*}
\det\left(
\begin{array}
[c]{c}%
A\leftarrow B\\
J\leftarrow K
\end{array}
\right)   &  =\sum_{\substack{P\subseteq\left\{  1,2,\ldots,n\right\}
;\\\left\vert P\right\vert =\left\vert J\right\vert }}\left(  -1\right)
^{\sum P+\sum J}\det\left(  \operatorname*{sub}\nolimits_{w\left(  P\right)
}^{w\left(  K\right)  }B\right)  \det\left(  \operatorname*{sub}%
\nolimits_{w\left(  \widetilde{P}\right)  }^{w\left(  \widetilde{J}\right)
}A\right) \\
&  =\sum_{\substack{R\subseteq\left\{  1,2,\ldots,n\right\}  ;\\\left\vert
R\right\vert =\left\vert J\right\vert }}\left(  -1\right)  ^{\sum R+\sum
J}\det\left(  \operatorname*{sub}\nolimits_{w\left(  R\right)  }^{w\left(
K\right)  }B\right)  \det\left(  \operatorname*{sub}\nolimits_{w\left(
\widetilde{R}\right)  }^{w\left(  \widetilde{J}\right)  }A\right) \\
&  \ \ \ \ \ \ \ \ \ \ \ \ \ \ \ \ \ \ \ \ \left(
\begin{array}
[c]{c}%
\text{here, we have renamed the}\\
\text{summation index }P\text{ as }R
\end{array}
\right)  .
\end{align*}
Also, Lemma \ref{lem.sol.det.syl-lem.1} (applied to $B$, $A$, $K$ and $J$
instead of $A$, $B$, $J$ and $K$) yields%
\begin{align*}
\det\left(
\begin{array}
[c]{c}%
B\leftarrow A\\
K\leftarrow J
\end{array}
\right)   &  =\sum_{\substack{P\subseteq\left\{  1,2,\ldots,n\right\}
;\\\left\vert P\right\vert =\left\vert K\right\vert }}\left(  -1\right)
^{\sum P+\sum K}\det\left(  \operatorname*{sub}\nolimits_{w\left(  P\right)
}^{w\left(  J\right)  }A\right)  \det\left(  \operatorname*{sub}%
\nolimits_{w\left(  \widetilde{P}\right)  }^{w\left(  \widetilde{K}\right)
}B\right) \\
&  \ \ \ \ \ \ \ \ \ \ \ \ \ \ \ \ \ \ \ \ \left(  \text{since }\left\vert
J\right\vert =\left\vert K\right\vert \text{ (because }\left\vert K\right\vert
=\left\vert J\right\vert \text{)}\right) \\
&  =\sum_{\substack{P\subseteq\left\{  1,2,\ldots,n\right\}  ;\\\left\vert
P\right\vert =\left\vert J\right\vert }}\left(  -1\right)  ^{\sum P+\sum
K}\det\left(  \operatorname*{sub}\nolimits_{w\left(  P\right)  }^{w\left(
J\right)  }A\right)  \det\left(  \operatorname*{sub}\nolimits_{w\left(
\widetilde{P}\right)  }^{w\left(  \widetilde{K}\right)  }B\right) \\
&  \ \ \ \ \ \ \ \ \ \ \ \ \ \ \ \ \ \ \ \ \left(  \text{since }\left\vert
K\right\vert =\left\vert J\right\vert \right)  .
\end{align*}
Multiplying these two equalities, we find%
\begin{align*}
&  \det\left(
\begin{array}
[c]{c}%
A\leftarrow B\\
J\leftarrow K
\end{array}
\right)  \cdot\det\left(
\begin{array}
[c]{c}%
B\leftarrow A\\
K\leftarrow J
\end{array}
\right) \\
&  =\left(  \sum_{\substack{R\subseteq\left\{  1,2,\ldots,n\right\}
;\\\left\vert R\right\vert =\left\vert J\right\vert }}\left(  -1\right)
^{\sum R+\sum J}\det\left(  \operatorname*{sub}\nolimits_{w\left(  R\right)
}^{w\left(  K\right)  }B\right)  \det\left(  \operatorname*{sub}%
\nolimits_{w\left(  \widetilde{R}\right)  }^{w\left(  \widetilde{J}\right)
}A\right)  \right) \\
&  \ \ \ \ \ \ \ \ \ \ \cdot\left(  \sum_{\substack{P\subseteq\left\{
1,2,\ldots,n\right\}  ;\\\left\vert P\right\vert =\left\vert J\right\vert
}}\left(  -1\right)  ^{\sum P+\sum K}\det\left(  \operatorname*{sub}%
\nolimits_{w\left(  P\right)  }^{w\left(  J\right)  }A\right)  \det\left(
\operatorname*{sub}\nolimits_{w\left(  \widetilde{P}\right)  }^{w\left(
\widetilde{K}\right)  }B\right)  \right) \\
&  =\sum_{\substack{R\subseteq\left\{  1,2,\ldots,n\right\}  ;\\\left\vert
R\right\vert =\left\vert J\right\vert }}\ \ \sum_{\substack{P\subseteq\left\{
1,2,\ldots,n\right\}  ;\\\left\vert P\right\vert =\left\vert J\right\vert
}}\left(  -1\right)  ^{\sum R+\sum J}\det\left(  \operatorname*{sub}%
\nolimits_{w\left(  R\right)  }^{w\left(  K\right)  }B\right)  \det\left(
\operatorname*{sub}\nolimits_{w\left(  \widetilde{R}\right)  }^{w\left(
\widetilde{J}\right)  }A\right) \\
&  \ \ \ \ \ \ \ \ \ \ \ \ \ \ \ \ \ \ \ \ \cdot\left(  -1\right)  ^{\sum
P+\sum K}\det\left(  \operatorname*{sub}\nolimits_{w\left(  P\right)
}^{w\left(  J\right)  }A\right)  \det\left(  \operatorname*{sub}%
\nolimits_{w\left(  \widetilde{P}\right)  }^{w\left(  \widetilde{K}\right)
}B\right) \\
&  =\sum_{\substack{R\subseteq\left\{  1,2,\ldots,n\right\}  ;\\\left\vert
R\right\vert =\left\vert J\right\vert }}\ \ \sum_{\substack{P\subseteq\left\{
1,2,\ldots,n\right\}  ;\\\left\vert P\right\vert =\left\vert J\right\vert
}}\underbrace{\left(  -1\right)  ^{\sum R+\sum J}\det\left(
\operatorname*{sub}\nolimits_{w\left(  P\right)  }^{w\left(  J\right)
}A\right)  \det\left(  \operatorname*{sub}\nolimits_{w\left(  \widetilde{R}%
\right)  }^{w\left(  \widetilde{J}\right)  }A\right)  }_{\substack{=\alpha
_{P,R}\\\text{(by (\ref{sol.det.syl-lem.beta=}))}}}\\
&  \ \ \ \ \ \ \ \ \ \ \ \ \ \ \ \ \ \ \ \ \cdot\left(  -1\right)  ^{\sum
P+\sum K}\det\left(  \operatorname*{sub}\nolimits_{w\left(  R\right)
}^{w\left(  K\right)  }B\right)  \det\left(  \operatorname*{sub}%
\nolimits_{w\left(  \widetilde{P}\right)  }^{w\left(  \widetilde{K}\right)
}B\right) \\
&  \ \ \ \ \ \ \ \ \ \ \ \ \ \ \ \ \ \ \ \ \ \ \ \ \ \ \ \ \ \ \left(
\text{here, we have permuted the factors in the product}\right) \\
&  =\sum_{\substack{R\subseteq\left\{  1,2,\ldots,n\right\}  ;\\\left\vert
R\right\vert =\left\vert J\right\vert }}\ \ \sum_{\substack{P\subseteq\left\{
1,2,\ldots,n\right\}  ;\\\left\vert P\right\vert =\left\vert J\right\vert
}}\alpha_{P,R}\cdot\left(  -1\right)  ^{\sum P+\sum K}\det\left(
\operatorname*{sub}\nolimits_{w\left(  R\right)  }^{w\left(  K\right)
}B\right)  \det\left(  \operatorname*{sub}\nolimits_{w\left(  \widetilde{P}%
\right)  }^{w\left(  \widetilde{K}\right)  }B\right)  .
\end{align*}

Forget that we fixed $K$. We thus have proved the equality%
\begin{align*}
&  \det\left(
\begin{array}
[c]{c}%
A\leftarrow B\\
J\leftarrow K
\end{array}
\right)  \cdot\det\left(
\begin{array}
[c]{c}%
B\leftarrow A\\
K\leftarrow J
\end{array}
\right) \\
&  =\sum_{\substack{R\subseteq\left\{  1,2,\ldots,n\right\}  ;\\\left\vert
R\right\vert =\left\vert J\right\vert }}\ \ \sum_{\substack{P\subseteq\left\{
1,2,\ldots,n\right\}  ;\\\left\vert P\right\vert =\left\vert J\right\vert
}}\alpha_{P,R}\cdot\left(  -1\right)  ^{\sum P+\sum K}\det\left(
\operatorname*{sub}\nolimits_{w\left(  R\right)  }^{w\left(  K\right)
}B\right)  \det\left(  \operatorname*{sub}\nolimits_{w\left(  \widetilde{P}%
\right)  }^{w\left(  \widetilde{K}\right)  }B\right)
\end{align*}
for any subset $K$ of $\left\{  1,2,\ldots,n\right\}  $ satisfying $\left\vert
K\right\vert =\left\vert J\right\vert $.

If we sum these equalities over all such subsets $K$, then we obtain
\begin{align}
&  \sum_{\substack{K\subseteq\left\{  1,2,\ldots,n\right\}  ;\\\left\vert
K\right\vert =\left\vert J\right\vert }}\det\left(
\begin{array}
[c]{c}%
A\leftarrow B\\
J\leftarrow K
\end{array}
\right)  \det\left(
\begin{array}
[c]{c}%
B\leftarrow A\\
K\leftarrow J
\end{array}
\right) \nonumber\\
&  =\underbrace{\sum_{\substack{K\subseteq\left\{  1,2,\ldots,n\right\}
;\\\left\vert K\right\vert =\left\vert J\right\vert }}\ \ \sum
_{\substack{R\subseteq\left\{  1,2,\ldots,n\right\}  ;\\\left\vert
R\right\vert =\left\vert J\right\vert }}\ \ \sum_{\substack{P\subseteq\left\{
1,2,\ldots,n\right\}  ;\\\left\vert P\right\vert =\left\vert J\right\vert }%
}}_{=\sum_{\substack{P\subseteq\left\{  1,2,\ldots,n\right\}  ;\\\left\vert
P\right\vert =\left\vert J\right\vert }}\ \ \sum_{\substack{R\subseteq\left\{
1,2,\ldots,n\right\}  ;\\\left\vert R\right\vert =\left\vert J\right\vert
}}\ \ \sum_{\substack{K\subseteq\left\{  1,2,\ldots,n\right\}  ;\\\left\vert
K\right\vert =\left\vert J\right\vert }}}\alpha_{P,R}\cdot\left(  -1\right)
^{\sum P+\sum K}\det\left(  \operatorname*{sub}\nolimits_{w\left(  R\right)
}^{w\left(  K\right)  }B\right)  \det\left(  \operatorname*{sub}%
\nolimits_{w\left(  \widetilde{P}\right)  }^{w\left(  \widetilde{K}\right)
}B\right) \nonumber\\
&  =\sum_{\substack{P\subseteq\left\{  1,2,\ldots,n\right\}  ;\\\left\vert
P\right\vert =\left\vert J\right\vert }}\ \ \sum_{\substack{R\subseteq\left\{
1,2,\ldots,n\right\}  ;\\\left\vert R\right\vert =\left\vert J\right\vert
}}\ \ \sum_{\substack{K\subseteq\left\{  1,2,\ldots,n\right\}  ;\\\left\vert
K\right\vert =\left\vert J\right\vert }}\alpha_{P,R}\cdot\left(  -1\right)
^{\sum P+\sum K}\det\left(  \operatorname*{sub}\nolimits_{w\left(  R\right)
}^{w\left(  K\right)  }B\right)  \det\left(  \operatorname*{sub}%
\nolimits_{w\left(  \widetilde{P}\right)  }^{w\left(  \widetilde{K}\right)
}B\right) \nonumber\\
&  =\sum_{\substack{P\subseteq\left\{  1,2,\ldots,n\right\}  ;\\\left\vert
P\right\vert =\left\vert J\right\vert }}\ \ \sum_{\substack{R\subseteq\left\{
1,2,\ldots,n\right\}  ;\\\left\vert R\right\vert =\left\vert J\right\vert
}}\ \ \sum_{\substack{Q\subseteq\left\{  1,2,\ldots,n\right\}  ;\\\left\vert
Q\right\vert =\left\vert J\right\vert }}\alpha_{P,R}\cdot\left(  -1\right)
^{\sum P+\sum Q}\det\left(  \operatorname*{sub}\nolimits_{w\left(  R\right)
}^{w\left(  Q\right)  }B\right)  \det\left(  \operatorname*{sub}%
\nolimits_{w\left(  \widetilde{P}\right)  }^{w\left(  \widetilde{Q}\right)
}B\right) \nonumber\\
&  \ \ \ \ \ \ \ \ \ \ \ \ \ \ \ \ \ \ \ \ \left(  \text{here, we have renamed
the summation index }K\text{ as }Q\right) \nonumber\\
&  =\sum_{\substack{P\subseteq\left\{  1,2,\ldots,n\right\}  ;\\\left\vert
P\right\vert =\left\vert J\right\vert }}\ \ \sum_{\substack{R\subseteq\left\{
1,2,\ldots,n\right\}  ;\\\left\vert R\right\vert =\left\vert J\right\vert
}}\alpha_{P,R}\sum_{\substack{Q\subseteq\left\{  1,2,\ldots,n\right\}
;\\\left\vert Q\right\vert =\left\vert J\right\vert }}\left(  -1\right)
^{\sum P+\sum Q}\det\left(  \operatorname*{sub}\nolimits_{w\left(  R\right)
}^{w\left(  Q\right)  }B\right)  \det\left(  \operatorname*{sub}%
\nolimits_{w\left(  \widetilde{P}\right)  }^{w\left(  \widetilde{Q}\right)
}B\right)  . \label{sol.det.syl-lem.5}%
\end{align}

We shall now analyze the inner sum on the right hand side of
(\ref{sol.det.syl-lem.5}).

First, we fix a subset $P$ of $\left\{  1,2,\ldots,n\right\}  $ satisfying
$\left\vert P\right\vert =\left\vert J\right\vert $. Then,%
\begin{align}
&  \sum_{\substack{R\subseteq\left\{  1,2,\ldots,n\right\}  ;\\\left\vert
R\right\vert =\left\vert J\right\vert }}\alpha_{P,R}\sum_{\substack{Q\subseteq
\left\{  1,2,\ldots,n\right\}  ;\\\left\vert Q\right\vert =\left\vert
J\right\vert }}\left(  -1\right)  ^{\sum P+\sum Q}\det\left(
\operatorname*{sub}\nolimits_{w\left(  R\right)  }^{w\left(  Q\right)
}B\right)  \det\left(  \operatorname*{sub}\nolimits_{w\left(  \widetilde{P}%
\right)  }^{w\left(  \widetilde{Q}\right)  }B\right) \nonumber\\
&  =\sum_{\substack{R\subseteq\left\{  1,2,\ldots,n\right\}  ;\\\left\vert
R\right\vert =\left\vert J\right\vert }}\alpha_{P,R}\sum_{\substack{Q\subseteq
\left\{  1,2,\ldots,n\right\}  ;\\\left\vert Q\right\vert =\left\vert
P\right\vert }}\left(  -1\right)  ^{\sum P+\sum Q}\det\left(
\operatorname*{sub}\nolimits_{w\left(  R\right)  }^{w\left(  Q\right)
}B\right)  \det\left(  \operatorname*{sub}\nolimits_{w\left(  \widetilde{P}%
\right)  }^{w\left(  \widetilde{Q}\right)  }B\right) \nonumber\\
&  \ \ \ \ \ \ \ \ \ \ \ \ \ \ \ \ \ \ \ \ \left(  \text{since }\left\vert
J\right\vert =\left\vert P\right\vert \right) \nonumber\\
&  =\alpha_{P,P}\underbrace{\sum_{\substack{Q\subseteq\left\{  1,2,\ldots
,n\right\}  ;\\\left\vert Q\right\vert =\left\vert P\right\vert }}\left(
-1\right)  ^{\sum P+\sum Q}\det\left(  \operatorname*{sub}\nolimits_{w\left(
P\right)  }^{w\left(  Q\right)  }B\right)  \det\left(  \operatorname*{sub}%
\nolimits_{w\left(  \widetilde{P}\right)  }^{w\left(  \widetilde{Q}\right)
}B\right)  }_{\substack{=\det B\\\text{(since Theorem
\ref{thm.det.laplace-multi} \textbf{(a)} (applied to }B\text{ instead of
}A\text{)}\\\text{yields }\det B=\sum_{\substack{Q\subseteq\left\{
1,2,\ldots,n\right\}  ;\\\left\vert Q\right\vert =\left\vert P\right\vert
}}\left(  -1\right)  ^{\sum P+\sum Q}\det\left(  \operatorname*{sub}%
\nolimits_{w\left(  P\right)  }^{w\left(  Q\right)  }B\right)  \det\left(
\operatorname*{sub}\nolimits_{w\left(  \widetilde{P}\right)  }^{w\left(
\widetilde{Q}\right)  }B\right)  \text{)}}}\nonumber\\
&  \ \ \ \ \ \ \ \ \ \ +\sum_{\substack{R\subseteq\left\{  1,2,\ldots
,n\right\}  ;\\\left\vert R\right\vert =\left\vert J\right\vert ;\\R\neq
P}}\alpha_{P,R}\underbrace{\sum_{\substack{Q\subseteq\left\{  1,2,\ldots
,n\right\}  ;\\\left\vert Q\right\vert =\left\vert P\right\vert }}\left(
-1\right)  ^{\sum P+\sum Q}\det\left(  \operatorname*{sub}\nolimits_{w\left(
R\right)  }^{w\left(  Q\right)  }B\right)  \det\left(  \operatorname*{sub}%
\nolimits_{w\left(  \widetilde{P}\right)  }^{w\left(  \widetilde{Q}\right)
}B\right)  }_{\substack{=0\\\text{(since Exercise
\ref{exe.det.laplace-multi.0} \textbf{(a)} (applied to }B\text{ instead of
}A\text{)}\\\text{yields }0=\sum_{\substack{Q\subseteq\left\{  1,2,\ldots
,n\right\}  ;\\\left\vert Q\right\vert =\left\vert P\right\vert }}\left(
-1\right)  ^{\sum P+\sum Q}\det\left(  \operatorname*{sub}\nolimits_{w\left(
R\right)  }^{w\left(  Q\right)  }B\right)  \det\left(  \operatorname*{sub}%
\nolimits_{w\left(  \widetilde{P}\right)  }^{w\left(  \widetilde{Q}\right)
}B\right)  \\\text{(because }\left\vert Q\right\vert =\left\vert P\right\vert
=\left\vert J\right\vert =\left\vert R\right\vert \text{ and }P\neq
R\text{))}}}\nonumber\\
&  \ \ \ \ \ \ \ \ \ \ \ \ \ \ \ \ \ \ \ \ \left(
\begin{array}
[c]{c}%
\text{here, we have split off the addend for }R=P\\
\text{from the sum (since }P\subseteq\left\{  1,2,\ldots,n\right\}  \text{ and
}\left\vert P\right\vert =\left\vert J\right\vert \text{)}%
\end{array}
\right) \nonumber\\
&  =\alpha_{P,P}\det B+\underbrace{\sum_{\substack{R\subseteq\left\{
1,2,\ldots,n\right\}  ;\\\left\vert R\right\vert =\left\vert J\right\vert
;\\R\neq P}}\alpha_{P,R}0}_{=0}=\alpha_{P,P}\det B. \label{sol.det.syl-lem.12}%
\end{align}

Forget that we fixed $P$. We thus have proved the equality
(\ref{sol.det.syl-lem.12}) for every subset $P$ of $\left\{  1,2,\ldots
,n\right\}  $ satisfying $\left\vert P\right\vert =\left\vert J\right\vert $.

Now, (\ref{sol.det.syl-lem.5}) becomes%
\begin{align*}
&  \sum_{\substack{K\subseteq\left\{  1,2,\ldots,n\right\}  ;\\\left\vert
K\right\vert =\left\vert J\right\vert }}\det\left(
\begin{array}
[c]{c}%
A\leftarrow B\\
J\leftarrow K
\end{array}
\right)  \det\left(
\begin{array}
[c]{c}%
B\leftarrow A\\
K\leftarrow J
\end{array}
\right) \\
&  =\sum_{\substack{P\subseteq\left\{  1,2,\ldots,n\right\}  ;\\\left\vert
P\right\vert =\left\vert J\right\vert }}\ \ \underbrace{\sum
_{\substack{R\subseteq\left\{  1,2,\ldots,n\right\}  ;\\\left\vert
R\right\vert =\left\vert J\right\vert }}\alpha_{P,R}\sum_{\substack{Q\subseteq
\left\{  1,2,\ldots,n\right\}  ;\\\left\vert Q\right\vert =\left\vert
J\right\vert }}\left(  -1\right)  ^{\sum P+\sum Q}\det\left(
\operatorname*{sub}\nolimits_{w\left(  R\right)  }^{w\left(  Q\right)
}B\right)  \det\left(  \operatorname*{sub}\nolimits_{w\left(  \widetilde{P}%
\right)  }^{w\left(  \widetilde{Q}\right)  }B\right)  }_{\substack{=\alpha
_{P,P}\det B\\\text{(by (\ref{sol.det.syl-lem.12}))}}}\\
&  =\sum_{\substack{P\subseteq\left\{  1,2,\ldots,n\right\}  ;\\\left\vert
P\right\vert =\left\vert J\right\vert }}\underbrace{\alpha_{P,P}%
}_{\substack{=\left(  -1\right)  ^{\sum P+\sum J}\det\left(
\operatorname*{sub}\nolimits_{w\left(  P\right)  }^{w\left(  J\right)
}A\right)  \det\left(  \operatorname*{sub}\nolimits_{w\left(  \widetilde{P}%
\right)  }^{w\left(  \widetilde{J}\right)  }A\right)  \\\text{(by the
definition of }\alpha_{P,P}\text{)}}}\det B\\
&  =\sum_{\substack{P\subseteq\left\{  1,2,\ldots,n\right\}  ;\\\left\vert
P\right\vert =\left\vert J\right\vert }}\left(  -1\right)  ^{\sum P+\sum
J}\det\left(  \operatorname*{sub}\nolimits_{w\left(  P\right)  }^{w\left(
J\right)  }A\right)  \det\left(  \operatorname*{sub}\nolimits_{w\left(
\widetilde{P}\right)  }^{w\left(  \widetilde{J}\right)  }A\right)  \det B\\
&  =\underbrace{\left(  \sum_{\substack{P\subseteq\left\{  1,2,\ldots
,n\right\}  ;\\\left\vert P\right\vert =\left\vert J\right\vert }}\left(
-1\right)  ^{\sum P+\sum J}\det\left(  \operatorname*{sub}\nolimits_{w\left(
P\right)  }^{w\left(  J\right)  }A\right)  \det\left(  \operatorname*{sub}%
\nolimits_{w\left(  \widetilde{P}\right)  }^{w\left(  \widetilde{J}\right)
}A\right)  \right)  }_{\substack{=\det A\\\text{(since Theorem
\ref{thm.det.laplace-multi} \textbf{(b)} (applied to }Q=J\text{)}%
\\\text{yields }\det A=\sum_{\substack{P\subseteq\left\{  1,2,\ldots
,n\right\}  ;\\\left\vert P\right\vert =\left\vert J\right\vert }}\left(
-1\right)  ^{\sum P+\sum J}\det\left(  \operatorname*{sub}\nolimits_{w\left(
P\right)  }^{w\left(  J\right)  }A\right)  \det\left(  \operatorname*{sub}%
\nolimits_{w\left(  \widetilde{P}\right)  }^{w\left(  \widetilde{J}\right)
}A\right)  \text{)}}}\cdot\det B\\
&  =\det A\cdot\det B.
\end{align*}
This solves Exercise \ref{exe.det.syl-lem}.
\end{proof}

\subsection{Solution to Exercise \ref{exe.det.muir-col-by-col}}

In this section, we shall use the notations introduced in Definition
\ref{def.sect.laplace.notations}. We begin with a lemma:

\begin{proof}
[Solution to Exercise \ref{exe.det.muir-col-by-col}.]For any subset $I$ of
$\left\{  1,2,\ldots,n\right\}  $, we let $\widetilde{I}$ denote the
complement $\left\{  1,2,\ldots,n\right\}  \setminus I$ of $I$. (For instance,
if $n=4$ and $I=\left\{  1,4\right\}  $, then $\widetilde{I}=\left\{
2,3\right\}  $.)

Let $J$ be a subset of $\left\{  1,2,\ldots,n\right\}  $ satisfying
$\left\vert J\right\vert =r$. Then, $\left\vert J\right\vert =\left\vert
J\right\vert $. Hence, Lemma \ref{lem.sol.det.syl-lem.1} (applied to $K=J$)
yields%
\begin{align}
&  \det\left(
\begin{array}
[c]{c}%
A\leftarrow B\\
J\leftarrow J
\end{array}
\right) \nonumber\\
&  =\sum_{\substack{P\subseteq\left\{  1,2,\ldots,n\right\}  ;\\\left\vert
P\right\vert =\left\vert J\right\vert }}\left(  -1\right)  ^{\sum P+\sum
J}\det\left(  \operatorname*{sub}\nolimits_{w\left(  P\right)  }^{w\left(
J\right)  }B\right)  \det\left(  \operatorname*{sub}\nolimits_{w\left(
\widetilde{P}\right)  }^{w\left(  \widetilde{J}\right)  }A\right) \nonumber\\
&  =\sum_{\substack{P\subseteq\left\{  1,2,\ldots,n\right\}  ;\\\left\vert
P\right\vert =r}}\left(  -1\right)  ^{\sum P+\sum J}\det\left(
\operatorname*{sub}\nolimits_{w\left(  P\right)  }^{w\left(  J\right)
}B\right)  \det\left(  \operatorname*{sub}\nolimits_{w\left(  \widetilde{P}%
\right)  }^{w\left(  \widetilde{J}\right)  }A\right)
\label{sol.det.muir-col-by-col.1}%
\end{align}
(since $\left\vert J\right\vert =r$). Furthermore, Lemma
\ref{lem.sol.det.syl-lem.1} (applied to $A^{T}$, $B^{T}$ and $J$ instead of
$A$, $B$ and $K$) yields
\begin{align}
&  \det\left(
\begin{array}
[c]{c}%
A^{T}\leftarrow B^{T}\\
J\leftarrow J
\end{array}
\right) \nonumber\\
&  =\sum_{\substack{P\subseteq\left\{  1,2,\ldots,n\right\}  ;\\\left\vert
P\right\vert =\left\vert J\right\vert }}\left(  -1\right)  ^{\sum P+\sum
J}\det\left(  \operatorname*{sub}\nolimits_{w\left(  P\right)  }^{w\left(
J\right)  }\left(  B^{T}\right)  \right)  \det\left(  \operatorname*{sub}%
\nolimits_{w\left(  \widetilde{P}\right)  }^{w\left(  \widetilde{J}\right)
}\left(  A^{T}\right)  \right) \nonumber\\
&  =\sum_{\substack{P\subseteq\left\{  1,2,\ldots,n\right\}  ;\\\left\vert
P\right\vert =r}}\left(  -1\right)  ^{\sum P+\sum J}\det\left(
\operatorname*{sub}\nolimits_{w\left(  P\right)  }^{w\left(  J\right)
}\left(  B^{T}\right)  \right)  \det\left(  \operatorname*{sub}%
\nolimits_{w\left(  \widetilde{P}\right)  }^{w\left(  \widetilde{J}\right)
}\left(  A^{T}\right)  \right)  \ \ \ \ \ \ \ \ \ \ \left(  \text{since
}\left\vert P\right\vert =r\right) \nonumber\\
&  =\sum_{\substack{Q\subseteq\left\{  1,2,\ldots,n\right\}  ;\\\left\vert
Q\right\vert =r}}\left(  -1\right)  ^{\sum Q+\sum J}\det\left(
\operatorname*{sub}\nolimits_{w\left(  Q\right)  }^{w\left(  J\right)
}\left(  B^{T}\right)  \right)  \det\left(  \operatorname*{sub}%
\nolimits_{w\left(  \widetilde{Q}\right)  }^{w\left(  \widetilde{J}\right)
}\left(  A^{T}\right)  \right)  \label{sol.det.muir-col-by-col.2}%
\end{align}
(here, we have renamed the summation index $P$ as $Q$).

Now, forget that we fixed $J$. We thus have proved the two equalities
(\ref{sol.det.muir-col-by-col.1}) and (\ref{sol.det.muir-col-by-col.2}) for
all subsets $J$ of $\left\{  1,2,\ldots,n\right\}  $ satisfying $\left\vert
J\right\vert =r$.

Now, let $J$ and $Q$ be two subsets of $\left\{  1,2,\ldots,n\right\}  $
satisfying $\left\vert J\right\vert =r$ and $\left\vert Q\right\vert =r$.
Then, $\left\vert J\right\vert =\left\vert Q\right\vert $. Hence, Corollary
\ref{cor.sol.det.laplace-multi.0.detsubAT} (applied to $n$, $Q$, $J$ and $B$
instead of $m$, $U$, $V$ and $A$) yields
\begin{equation}
\det\left(  \operatorname*{sub}\nolimits_{w\left(  Q\right)  }^{w\left(
J\right)  }\left(  B^{T}\right)  \right)  =\det\left(  \operatorname*{sub}%
\nolimits_{w\left(  J\right)  }^{w\left(  Q\right)  }B\right)  .
\label{sol.det.muir-col-by-col.3}%
\end{equation}
Furthermore, the definition of $\widetilde{J}$ yields $\widetilde{J}=\left\{
1,2,\ldots,n\right\}  \setminus J$. Thus,%
\begin{align*}
\left\vert \widetilde{J}\right\vert  &  =\left\vert \left\{  1,2,\ldots
,n\right\}  \setminus J\right\vert =\left\vert \left\{  1,2,\ldots,n\right\}
\right\vert -\left\vert J\right\vert \ \ \ \ \ \ \ \ \ \ \left(  \text{since
}J\text{ is a subset of }\left\{  1,2,\ldots,n\right\}  \right) \\
&  =n-\left\vert J\right\vert \ \ \ \ \ \ \ \ \ \ \left(  \text{since
}\left\vert \left\{  1,2,\ldots,n\right\}  \right\vert =n\right)  .
\end{align*}
Similarly, $\left\vert \widetilde{Q}\right\vert =n-\left\vert Q\right\vert $.
Now, $\left\vert \widetilde{J}\right\vert =n-\underbrace{\left\vert
J\right\vert }_{=\left\vert Q\right\vert }=n-\left\vert Q\right\vert
=\left\vert \widetilde{Q}\right\vert $. Hence, Corollary
\ref{cor.sol.det.laplace-multi.0.detsubAT} (applied to $n$, $\widetilde{Q}$
and $\widetilde{J}$ instead of $m$, $U$ and $V$) yields
\begin{equation}
\det\left(  \operatorname*{sub}\nolimits_{w\left(  \widetilde{Q}\right)
}^{w\left(  \widetilde{J}\right)  }\left(  A^{T}\right)  \right)  =\det\left(
\operatorname*{sub}\nolimits_{w\left(  \widetilde{J}\right)  }^{w\left(
\widetilde{Q}\right)  }A\right)  . \label{sol.det.muir-col-by-col.4}%
\end{equation}

Now, forget that we fixed $J$ and $Q$. We thus have proved the two equalities
(\ref{sol.det.muir-col-by-col.3}) and (\ref{sol.det.muir-col-by-col.4})
whenever $J$ and $Q$ are two subsets of $\left\{  1,2,\ldots,n\right\}  $
satisfying $\left\vert J\right\vert =r$ and $\left\vert Q\right\vert =r$.

Now,
\begin{align*}
&  \sum_{\substack{J\subseteq\left\{  1,2,\ldots,n\right\}  ;\\\left\vert
J\right\vert =r}}\underbrace{\det\left(
\begin{array}
[c]{c}%
A^{T}\leftarrow B^{T}\\
J\leftarrow J
\end{array}
\right)  }_{\substack{=\sum_{\substack{Q\subseteq\left\{  1,2,\ldots
,n\right\}  ;\\\left\vert Q\right\vert =r}}\left(  -1\right)  ^{\sum Q+\sum
J}\det\left(  \operatorname*{sub}\nolimits_{w\left(  Q\right)  }^{w\left(
J\right)  }\left(  B^{T}\right)  \right)  \det\left(  \operatorname*{sub}%
\nolimits_{w\left(  \widetilde{Q}\right)  }^{w\left(  \widetilde{J}\right)
}\left(  A^{T}\right)  \right)  \\\text{(by (\ref{sol.det.muir-col-by-col.2}%
))}}}\\
&  =\underbrace{\sum_{\substack{J\subseteq\left\{  1,2,\ldots,n\right\}
;\\\left\vert J\right\vert =r}}\ \ \sum_{\substack{Q\subseteq\left\{
1,2,\ldots,n\right\}  ;\\\left\vert Q\right\vert =r}}}_{=\sum
_{\substack{Q\subseteq\left\{  1,2,\ldots,n\right\}  ;\\\left\vert
Q\right\vert =r}}\ \ \sum_{\substack{J\subseteq\left\{  1,2,\ldots,n\right\}
;\\\left\vert J\right\vert =r}}}\underbrace{\left(  -1\right)  ^{\sum Q+\sum
J}}_{=\left(  -1\right)  ^{\sum J+\sum Q}}\underbrace{\det\left(
\operatorname*{sub}\nolimits_{w\left(  Q\right)  }^{w\left(  J\right)
}\left(  B^{T}\right)  \right)  }_{\substack{=\det\left(  \operatorname*{sub}%
\nolimits_{w\left(  J\right)  }^{w\left(  Q\right)  }B\right)  \\\text{(by
(\ref{sol.det.muir-col-by-col.3}))}}}\underbrace{\det\left(
\operatorname*{sub}\nolimits_{w\left(  \widetilde{Q}\right)  }^{w\left(
\widetilde{J}\right)  }\left(  A^{T}\right)  \right)  }_{\substack{=\det
\left(  \operatorname*{sub}\nolimits_{w\left(  \widetilde{J}\right)
}^{w\left(  \widetilde{Q}\right)  }A\right)  \\\text{(by
(\ref{sol.det.muir-col-by-col.4}))}}}\\
&  =\sum_{\substack{Q\subseteq\left\{  1,2,\ldots,n\right\}  ;\\\left\vert
Q\right\vert =r}}\ \ \sum_{\substack{J\subseteq\left\{  1,2,\ldots,n\right\}
;\\\left\vert J\right\vert =r}}\left(  -1\right)  ^{\sum J+\sum Q}\det\left(
\operatorname*{sub}\nolimits_{w\left(  J\right)  }^{w\left(  Q\right)
}B\right)  \det\left(  \operatorname*{sub}\nolimits_{w\left(  \widetilde{J}%
\right)  }^{w\left(  \widetilde{Q}\right)  }A\right) \\
&  =\sum_{\substack{Q\subseteq\left\{  1,2,\ldots,n\right\}  ;\\\left\vert
Q\right\vert =r}}\ \ \sum_{\substack{P\subseteq\left\{  1,2,\ldots,n\right\}
;\\\left\vert P\right\vert =r}}\left(  -1\right)  ^{\sum P+\sum Q}\det\left(
\operatorname*{sub}\nolimits_{w\left(  P\right)  }^{w\left(  Q\right)
}B\right)  \det\left(  \operatorname*{sub}\nolimits_{w\left(  \widetilde{P}%
\right)  }^{w\left(  \widetilde{Q}\right)  }A\right) \\
&  \ \ \ \ \ \ \ \ \ \ \ \ \ \ \ \ \ \ \ \ \left(
\begin{array}
[c]{c}%
\text{here, we have renamed the summation index }J\\
\text{as }P\text{ in the inner sum}%
\end{array}
\right) \\
&  =\sum_{\substack{J\subseteq\left\{  1,2,\ldots,n\right\}  ;\\\left\vert
J\right\vert =r}}\ \ \underbrace{\sum_{\substack{P\subseteq\left\{
1,2,\ldots,n\right\}  ;\\\left\vert P\right\vert =r}}\left(  -1\right)  ^{\sum
P+\sum J}\det\left(  \operatorname*{sub}\nolimits_{w\left(  P\right)
}^{w\left(  J\right)  }B\right)  \det\left(  \operatorname*{sub}%
\nolimits_{w\left(  \widetilde{P}\right)  }^{w\left(  \widetilde{J}\right)
}A\right)  }_{\substack{=\det\left(
\begin{array}
[c]{c}%
A\leftarrow B\\
J\leftarrow J
\end{array}
\right)  \\\text{(by (\ref{sol.det.muir-col-by-col.1}))}}}\\
&  \ \ \ \ \ \ \ \ \ \ \ \ \ \ \ \ \ \ \ \ \left(
\begin{array}
[c]{c}%
\text{here, we have renamed the summation index }Q\\
\text{as }J\text{ in the outer sum}%
\end{array}
\right) \\
&  =\sum_{\substack{J\subseteq\left\{  1,2,\ldots,n\right\}  ;\\\left\vert
J\right\vert =r}}\det\left(
\begin{array}
[c]{c}%
A\leftarrow B\\
J\leftarrow J
\end{array}
\right)  .
\end{align*}
This solves Exercise \ref{exe.det.muir-col-by-col}.
\end{proof}

\subsection{Solution to Exercise \ref{exe.det.AB=0-proport}}

\begin{proof}
[Solution to Exercise \ref{exe.det.AB=0-proport}.]We shall use the notations
introduced in Definition \ref{def.sect.laplace.notations} and in Definition
\ref{def.iverson}. From $k\in\left\{  0,1,\ldots,n\right\}  $, we obtain
$k\leq n$, so that $n-k\in\mathbb{N}$.

\begin{vershort}
The definition of the $n\times n$ identity matrix $I_{n}$ says that
$I_{n}=\left(  \delta_{i,j}\right)  _{1\leq i\leq n,\ 1\leq j\leq n}$, where
$\delta_{i,j}$ is defined to be $%
\begin{cases}
1, & \text{if }i=j;\\
0, & \text{if }i\neq j
\end{cases}
$. Using the Iverson bracket notation (see Definition \ref{def.rowscols}), we
can rewrite the number $\delta_{i,j}$ in this definition as $\left[
i=j\right]  $ (since both $\delta_{i,j}$ and $\left[  i=j\right]  $ are
defined to be $%
\begin{cases}
1, & \text{if }i=j;\\
0, & \text{if }i\neq j
\end{cases}
$). Hence, the definition of $I_{n}$ becomes%
\[
I_{n}=\left(  \left[  i=j\right]  \right)  _{1\leq i\leq n,\ 1\leq j\leq n}.
\]

\end{vershort}

\begin{verlong}
It is easy to see that%
\begin{equation}
I_{n}=\left(  \left[  i=j\right]  \right)  _{1\leq i\leq n,\ 1\leq j\leq n}
\label{sol.det.AB=0-proport.I=}%
\end{equation}
\footnote{\textit{Proof of (\ref{sol.det.AB=0-proport.I=}):} The definition of
the identity matrix $I_{n}$ yields $I_{n}=\left(  \delta_{i,j}\right)  _{1\leq
i\leq n,\ 1\leq j\leq n}$, where $\delta_{i,j}$ is defined to be $%
\begin{cases}
1, & \text{if }i=j;\\
0, & \text{if }i\neq j
\end{cases}
$. We can rewrite this as $I_{n}=\left(  \left[  i=j\right]  \right)  _{1\leq
i\leq n,\ 1\leq j\leq n}$, because every two elements $i\in\left\{
1,2,\ldots,n\right\}  $ and $j\in\left\{  1,2,\ldots,n\right\}  $ satisfy%
\begin{align*}
\delta_{i,j}  &  =%
\begin{cases}
1, & \text{if }i=j;\\
0, & \text{if }i\neq j
\end{cases}
\ \ \ \ \ \ \ \ \ \ \left(  \text{by the definition of }\delta_{i,j}\right) \\
&  =%
\begin{cases}
1, & \text{if }i=j\text{ is true};\\
0, & \text{if }i=j\text{ is false}%
\end{cases}
\ \ \ \ \ \ \ \ \ \ \left(  \text{since \textquotedblleft}i\neq
j\text{\textquotedblright\ means \textquotedblleft}i=j\text{ is
false\textquotedblright}\right) \\
&  =\left[  i=j\right]  \ \ \ \ \ \ \ \ \ \ \left(  \text{since }\left[
i=j\right]  \text{ is defined to be }%
\begin{cases}
1, & \text{if }i=j\text{ is true};\\
0, & \text{if }i=j\text{ is false}%
\end{cases}
\right)  .
\end{align*}
Thus, (\ref{sol.det.AB=0-proport.I=}) is proved.}.
\end{verlong}

The same argument (applied to $k$ instead of $n$) yields%
\begin{equation}
I_{k}=\left(  \left[  i=j\right]  \right)  _{1\leq i\leq k,\ 1\leq j\leq k}.
\label{sol.det.AB=0-proport.Ik=}%
\end{equation}
The same argument (applied to $n-k$ instead of $k$) yields%
\begin{equation}
I_{n-k}=\left(  \left[  i=j\right]  \right)  _{1\leq i\leq n-k,\ 1\leq j\leq
n-k} \label{sol.det.AB=0-proport.In-k=}%
\end{equation}
(since $n-k\in\mathbb{N}$).

\begin{vershort}
From Example \ref{exa.det.In}, we know that $\det\left(  I_{n}\right)  =1$.
Similarly, $\det\left(  I_{k}\right)  =1$ and $\det\left(  I_{n-k}\right)  =1$.
\end{vershort}

\begin{verlong}
Also, $\det\left(  I_{n}\right)  =1$ (by Example \ref{exa.det.In}). The same
argument (applied to $k$ instead of $n$) yields $\det\left(  I_{k}\right)
=1$. The same argument (applied to $n-k$ instead of $k$) yields $\det\left(
I_{n-k}\right)  =1$.
\end{verlong}

Write the $k\times n$-matrix $A$ as $A=\left(  a_{i,j}\right)  _{1\leq i\leq
k,\ 1\leq j\leq n}$. Write the $n\times\left(  n-k\right)  $-matrix $B$ as
$B=\left(  b_{i,j}\right)  _{1\leq i\leq n,\ 1\leq j\leq n-k}$. By assumption,
we have $AB=0_{k\times\left(  n-k\right)  }=\left(  0\right)  _{1\leq i\leq
k,\ 1\leq j\leq n-k}$. Hence,%
\[
\left(  0\right)  _{1\leq i\leq k,\ 1\leq j\leq n-k}=AB=\left(  \sum_{m=1}%
^{n}a_{i,m}b_{m,j}\right)  _{1\leq i\leq k,\ 1\leq j\leq n-k}%
\]
(by the definition of the product of two matrices, since $A=\left(
a_{i,j}\right)  _{1\leq i\leq k,\ 1\leq j\leq n}$ and $B=\left(
b_{i,j}\right)  _{1\leq i\leq n,\ 1\leq j\leq n-k}$). In other words, for any
$i\in\left\{  1,2,\ldots,k\right\}  $ and $j\in\left\{  1,2,\ldots
,n-k\right\}  $, we have%
\begin{equation}
0=\sum_{m=1}^{n}a_{i,m}b_{m,j}. \label{sol.det.AB=0-proport.0=sumab}%
\end{equation}

Proposition \ref{prop.sect.laplace.notations.w(I)} (applied to $I=P$) yields
that $w\left(  P\right)  $ is a $\left\vert P\right\vert $-tuple of elements
of $P$. In other words, $w\left(  P\right)  $ is a $k$-tuple of elements of
$P$ (since $\left\vert P\right\vert =k$). Write this $k$-tuple $w\left(
P\right)  $ as $w\left(  P\right)  =\left(  p_{1},p_{2},\ldots,p_{k}\right)
$. From $w\left(  P\right)  =\left(  p_{1},p_{2},\ldots,p_{k}\right)  $, we
obtain%
\begin{align}
\operatorname*{cols}\nolimits_{w\left(  P\right)  }A  &  =\operatorname*{cols}%
\nolimits_{\left(  p_{1},p_{2},\ldots,p_{k}\right)  }A\nonumber\\
&  =\operatorname*{cols}\nolimits_{p_{1},p_{2},\ldots,p_{k}}%
A\ \ \ \ \ \ \ \ \ \ \left(  \text{by the definition of }\operatorname*{cols}%
\nolimits_{\left(  p_{1},p_{2},\ldots,p_{k}\right)  }A\right) \nonumber\\
&  =\left(  a_{i,p_{y}}\right)  _{1\leq i\leq k,\ 1\leq y\leq k}
\label{sol.det.AB=0-proport.colswPA=}%
\end{align}
(by Definition \ref{def.rowscols} \textbf{(b)}, since $A=\left(
a_{i,j}\right)  _{1\leq i\leq k,\ 1\leq j\leq n}$).

\begin{vershort}
Proposition \ref{prop.sect.laplace.notations.w(I)} (applied to $I=Q$) yields
that $w\left(  Q\right)  $ is a $\left\vert Q\right\vert $-tuple of elements
of $Q$. In other words, $w\left(  Q\right)  $ is a $k$-tuple of elements of
$Q$ (since $\left\vert Q\right\vert =k$). Write this $k$-tuple $w\left(
Q\right)  $ as $w\left(  Q\right)  =\left(  q_{1},q_{2},\ldots,q_{k}\right)
$. Just as we proved (\ref{sol.det.AB=0-proport.colswPA=}), we can now see
that%
\begin{equation}
\operatorname*{cols}\nolimits_{w\left(  Q\right)  }A=\left(  a_{i,q_{y}%
}\right)  _{1\leq i\leq k,\ 1\leq y\leq k}.
\label{sol.det.AB=0-proport.colswQA=.short}%
\end{equation}

\end{vershort}

\begin{verlong}
Proposition \ref{prop.sect.laplace.notations.w(I)} (applied to $I=Q$) yields
that $w\left(  Q\right)  $ is a $\left\vert Q\right\vert $-tuple of elements
of $Q$. In other words, $w\left(  Q\right)  $ is a $k$-tuple of elements of
$Q$ (since $\left\vert Q\right\vert =k$). Write this $k$-tuple $w\left(
Q\right)  $ as $w\left(  Q\right)  =\left(  q_{1},q_{2},\ldots,q_{k}\right)
$. From $w\left(  Q\right)  =\left(  q_{1},q_{2},\ldots,q_{k}\right)  $, we
obtain%
\begin{align}
\operatorname*{cols}\nolimits_{w\left(  Q\right)  }A  &  =\operatorname*{cols}%
\nolimits_{\left(  q_{1},q_{2},\ldots,q_{k}\right)  }A\nonumber\\
&  =\operatorname*{cols}\nolimits_{q_{1},q_{2},\ldots,q_{k}}%
A\ \ \ \ \ \ \ \ \ \ \left(  \text{by the definition of }\operatorname*{cols}%
\nolimits_{\left(  q_{1},q_{2},\ldots,q_{k}\right)  }A\right) \nonumber\\
&  =\left(  a_{i,q_{y}}\right)  _{1\leq i\leq k,\ 1\leq y\leq k}
\label{sol.det.AB=0-proport.colswQA=}%
\end{align}
(by Definition \ref{def.rowscols} \textbf{(b)}, since $A=\left(
a_{i,j}\right)  _{1\leq i\leq k,\ 1\leq j\leq n}$).
\end{verlong}

\begin{vershort}
The set $P$ is a $k$-element subset of $\left\{  1,2,\ldots,n\right\}  $
(since $\left\vert P\right\vert =k$). Hence, its complement $\widetilde{P}$ is
an $\left(  n-k\right)  $-element subset of $\left\{  1,2,\ldots,n\right\}  $.
Thus, $\left\vert \widetilde{P}\right\vert =n-k$. Similarly, $\left\vert
\widetilde{Q}\right\vert =n-k$.
\end{vershort}

\begin{verlong}
The definition of $\widetilde{P}$ yields $\widetilde{P}=\left\{
1,2,\ldots,n\right\}  \setminus P$. Thus,%
\begin{align*}
\left\vert \widetilde{P}\right\vert  &  =\left\vert \left\{  1,2,\ldots
,n\right\}  \setminus P\right\vert =\underbrace{\left\vert \left\{
1,2,\ldots,n\right\}  \right\vert }_{=n}-\underbrace{\left\vert P\right\vert
}_{=k}\ \ \ \ \ \ \ \ \ \ \left(  \text{since }P\text{ is a subset of
}\left\{  1,2,\ldots,n\right\}  \right) \\
&  =n-k.
\end{align*}
Similarly, $\left\vert \widetilde{Q}\right\vert =n-k$.
\end{verlong}

Proposition \ref{prop.sect.laplace.notations.w(I)} (applied to
$I=\widetilde{P}$) yields that $w\left(  \widetilde{P}\right)  $ is a
$\left\vert \widetilde{P}\right\vert $-tuple of elements of $\widetilde{P}$.
In other words, $w\left(  \widetilde{P}\right)  $ is an $\left(  n-k\right)
$-tuple of elements of $\widetilde{P}$ (since $\left\vert \widetilde{P}%
\right\vert =n-k$). Write this $\left(  n-k\right)  $-tuple $w\left(
\widetilde{P}\right)  $ as $w\left(  \widetilde{P}\right)  =\left(
\widetilde{p}_{1},\widetilde{p}_{2},\ldots,\widetilde{p}_{n-k}\right)  $. From
$w\left(  \widetilde{P}\right)  =\left(  \widetilde{p}_{1},\widetilde{p}%
_{2},\ldots,\widetilde{p}_{n-k}\right)  $, we obtain%
\begin{align}
\operatorname*{rows}\nolimits_{w\left(  \widetilde{P}\right)  }B  &
=\operatorname*{rows}\nolimits_{\left(  \widetilde{p}_{1},\widetilde{p}%
_{2},\ldots,\widetilde{p}_{n-k}\right)  }B\nonumber\\
&  =\operatorname*{rows}\nolimits_{\widetilde{p}_{1},\widetilde{p}_{2}%
,\ldots,\widetilde{p}_{n-k}}B\ \ \ \ \ \ \ \ \ \ \left(  \text{by the
definition of }\operatorname*{rows}\nolimits_{\left(  \widetilde{p}%
_{1},\widetilde{p}_{2},\ldots,\widetilde{p}_{n-k}\right)  }B\right)
\nonumber\\
&  =\left(  b_{\widetilde{p}_{x},j}\right)  _{1\leq x\leq n-k,\ 1\leq j\leq
n-k} \label{sol.det.AB=0-proport.rowswPB=}%
\end{align}
(by Definition \ref{def.rowscols} \textbf{(a)}, since $B=\left(
b_{i,j}\right)  _{1\leq i\leq n,\ 1\leq j\leq n-k}$).

\begin{vershort}
Proposition \ref{prop.sect.laplace.notations.w(I)} (applied to
$I=\widetilde{Q}$) yields that $w\left(  \widetilde{Q}\right)  $ is a
$\left\vert \widetilde{Q}\right\vert $-tuple of elements of $\widetilde{Q}$.
In other words, $w\left(  \widetilde{Q}\right)  $ is an $\left(  n-k\right)
$-tuple of elements of $\widetilde{Q}$ (since $\left\vert \widetilde{Q}%
\right\vert =n-k$). Write this $\left(  n-k\right)  $-tuple $w\left(
\widetilde{Q}\right)  $ as $w\left(  \widetilde{Q}\right)  =\left(
\widetilde{q}_{1},\widetilde{q}_{2},\ldots,\widetilde{q}_{n-k}\right)  $. Just
as we proved (\ref{sol.det.AB=0-proport.rowswPB=}), we can now see that%
\begin{equation}
\operatorname*{rows}\nolimits_{w\left(  \widetilde{Q}\right)  }B=\left(
b_{\widetilde{q}_{x},j}\right)  _{1\leq x\leq n-k,\ 1\leq j\leq n-k}.
\label{sol.det.AB=0-proport.rowswQB=.short}%
\end{equation}

\end{vershort}

\begin{verlong}
Proposition \ref{prop.sect.laplace.notations.w(I)} (applied to
$I=\widetilde{Q}$) yields that $w\left(  \widetilde{Q}\right)  $ is a
$\left\vert \widetilde{Q}\right\vert $-tuple of elements of $\widetilde{Q}$.
In other words, $w\left(  \widetilde{Q}\right)  $ is an $\left(  n-k\right)
$-tuple of elements of $\widetilde{Q}$ (since $\left\vert \widetilde{Q}%
\right\vert =n-k$). Write this $\left(  n-k\right)  $-tuple $w\left(
\widetilde{Q}\right)  $ as $w\left(  \widetilde{Q}\right)  =\left(
\widetilde{q}_{1},\widetilde{q}_{2},\ldots,\widetilde{q}_{n-k}\right)  $. From
$w\left(  \widetilde{Q}\right)  =\left(  \widetilde{q}_{1},\widetilde{q}%
_{2},\ldots,\widetilde{q}_{n-k}\right)  $, we obtain%
\begin{align}
\operatorname*{rows}\nolimits_{w\left(  \widetilde{Q}\right)  }B  &
=\operatorname*{rows}\nolimits_{\left(  \widetilde{q}_{1},\widetilde{q}%
_{2},\ldots,\widetilde{q}_{n-k}\right)  }B\nonumber\\
&  =\operatorname*{rows}\nolimits_{\widetilde{q}_{1},\widetilde{q}_{2}%
,\ldots,\widetilde{q}_{n-k}}B\ \ \ \ \ \ \ \ \ \ \left(  \text{by the
definition of }\operatorname*{rows}\nolimits_{\left(  \widetilde{q}%
_{1},\widetilde{q}_{2},\ldots,\widetilde{q}_{n-k}\right)  }B\right)
\nonumber\\
&  =\left(  b_{\widetilde{q}_{x},j}\right)  _{1\leq x\leq n-k,\ 1\leq j\leq
n-k} \label{sol.det.AB=0-proport.rowswQB=}%
\end{align}
(by Definition \ref{def.rowscols} \textbf{(a)}, since $B=\left(
b_{i,j}\right)  _{1\leq i\leq n,\ 1\leq j\leq n-k}$).
\end{verlong}

\begin{vershort}
For any $i\in\left\{  1,2,\ldots,n\right\}  $ and $j\in\left\{  1,2,\ldots
,n\right\}  $, define an element $c_{i,j}\in\mathbb{K}$ by
\[
c_{i,j}=%
\begin{cases}
a_{i,j}, & \text{if }i\leq k;\\
\left[  j=\widetilde{p}_{i-k}\right]  , & \text{if }i>k
\end{cases}
\]
(where we are using the Iverson bracket notation). Define an $n\times
n$-matrix $C\in\mathbb{K}^{n\times n}$ by%
\[
C=\left(  c_{i,j}\right)  _{1\leq i\leq n,\ 1\leq j\leq n}.
\]

\end{vershort}

\begin{verlong}
For any $i\in\left\{  1,2,\ldots,n\right\}  $ and $j\in\left\{  1,2,\ldots
,n\right\}  $, define an element $c_{i,j}\in\mathbb{K}$ by
\[
c_{i,j}=%
\begin{cases}
a_{i,j}, & \text{if }i\leq k;\\
\left[  j=\widetilde{p}_{i-k}\right]  , & \text{if }i>k
\end{cases}
\]
(where we are using the Iverson bracket notation)\footnote{Let us check that
this is well-defined. Indeed:
\par
\begin{itemize}
\item If $i\in\left\{  1,2,\ldots,n\right\}  $ and $j\in\left\{
1,2,\ldots,n\right\}  $ satisfy $i\leq k$, then $a_{i,j}$ is well-defined
(since $i\in\left\{  1,2,\ldots,k\right\}  $ (because $i\in\left\{
1,2,\ldots,n\right\}  $ and $i\leq k$) and $j\in\left\{  1,2,\ldots,n\right\}
$).
\par
\item If $i\in\left\{  1,2,\ldots,n\right\}  $ and $j\in\left\{
1,2,\ldots,n\right\}  $ satisfy $i>k$, then $\widetilde{p}_{i-k}$ is
well-defined (since $i\in\left\{  k+1,k+2,\ldots,n\right\}  $ (because
$i\in\left\{  1,2,\ldots,n\right\}  $ and $i>k$) and therefore $i-k\in\left\{
1,2,\ldots,n-k\right\}  $) and therefore $\left[  j=\widetilde{p}%
_{i-k}\right]  $ is well-defined.
\end{itemize}
\par
Thus, $%
\begin{cases}
a_{i,j}, & \text{if }i\leq k;\\
\left[  j=\widetilde{p}_{i-k}\right]  , & \text{if }i>k
\end{cases}
$ is always well-defined.}. Define an $n\times n$-matrix $C\in\mathbb{K}%
^{n\times n}$ by%
\[
C=\left(  c_{i,j}\right)  _{1\leq i\leq n,\ 1\leq j\leq n}.
\]

\end{verlong}

\begin{vershort}
For any $i\in\left\{  1,2,\ldots,n\right\}  $ and $j\in\left\{  1,2,\ldots
,n\right\}  $, define an element $d_{i,j}\in\mathbb{K}$ by
\[
d_{i,j}=%
\begin{cases}
\left[  i=q_{j}\right]  , & \text{if }j\leq k;\\
b_{i,j-k}, & \text{if }j>k
\end{cases}
\]
(where we are using the Iverson bracket notation). Define an $n\times
n$-matrix $D\in\mathbb{K}^{n\times n}$ by%
\[
D=\left(  d_{i,j}\right)  _{1\leq i\leq n,\ 1\leq j\leq n}.
\]

\end{vershort}

\begin{verlong}
For any $i\in\left\{  1,2,\ldots,n\right\}  $ and $j\in\left\{  1,2,\ldots
,n\right\}  $, define an element $d_{i,j}\in\mathbb{K}$ by
\[
d_{i,j}=%
\begin{cases}
\left[  i=q_{j}\right]  , & \text{if }j\leq k;\\
b_{i,j-k}, & \text{if }j>k
\end{cases}
\]
(where we are using the Iverson bracket notation)\footnote{Let us check that
this is well-defined. Indeed:
\par
\begin{itemize}
\item If $i\in\left\{  1,2,\ldots,n\right\}  $ and $j\in\left\{
1,2,\ldots,n\right\}  $ satisfy $j\leq k$, then $q_{j}$ is well-defined (since
$j\in\left\{  1,2,\ldots,k\right\}  $ (because $j\in\left\{  1,2,\ldots
,n\right\}  $ and $j\leq k$)) and therefore $\left[  i=q_{j}\right]  $ is
well-defined.
\par
\item If $i\in\left\{  1,2,\ldots,n\right\}  $ and $j\in\left\{
1,2,\ldots,n\right\}  $ satisfy $j>k$, then $b_{i,j-k}$ is well-defined (since
$i\in\left\{  1,2,\ldots,n\right\}  $ and $j-k\in\left\{  1,2,\ldots
,n-k\right\}  $ (because we can combine $j\in\left\{  1,2,\ldots,n\right\}  $
with $j>k$ and thus obtain $j\in\left\{  k+1,k+2,\ldots,n\right\}  $)).
\end{itemize}
\par
Thus, $%
\begin{cases}
\left[  i=q_{j}\right]  , & \text{if }j\leq k;\\
b_{i,j-k}, & \text{if }j>k
\end{cases}
$ is always well-defined.}. Define an $n\times n$-matrix $D\in\mathbb{K}%
^{n\times n}$ by%
\[
D=\left(  d_{i,j}\right)  _{1\leq i\leq n,\ 1\leq j\leq n}.
\]

\end{verlong}

Theorem \ref{thm.det(AB)} (applied to $C$ and $D$ instead of $A$ and $B$)
yields%
\begin{equation}
\det\left(  CD\right)  =\det C\cdot\det D. \label{sol.det.AB=0-proport.detCD=}%
\end{equation}

Write the $n\times n$-matrix $CD\in\mathbb{K}^{n\times n}$ as $CD=\left(
g_{i,j}\right)  _{1\leq i\leq n,\ 1\leq j\leq n}$. Thus,%
\[
\left(  g_{i,j}\right)  _{1\leq i\leq n,\ 1\leq j\leq n}=CD=\left(  \sum
_{m=1}^{n}c_{i,m}d_{m,j}\right)  _{1\leq i\leq n,\ 1\leq j\leq n}%
\]
(by the definition of the product of two matrices, since $C=\left(
c_{i,j}\right)  _{1\leq i\leq n,\ 1\leq j\leq n}$ and $D=\left(
d_{i,j}\right)  _{1\leq i\leq n,\ 1\leq j\leq n}$). In other words, for any
$i\in\left\{  1,2,\ldots,n\right\}  $ and any $j\in\left\{  1,2,\ldots
,n\right\}  $, we have%
\begin{equation}
g_{i,j}=\sum_{m=1}^{n}c_{i,m}d_{m,j}. \label{sol.det.AB=0-proport.gij=}%
\end{equation}

\begin{vershort}
Let $K$ be the set $\left\{  1,2,\ldots,k\right\}  $. This is a subset of
$\left\{  1,2,\ldots,n\right\}  $, since $k\in\left\{  0,1,\ldots,n\right\}
$. Its size is $\left\vert K\right\vert =k$. Its complement is $\widetilde{K}%
=\left\{  k+1,k+2,\ldots,n\right\}  $. Hence, $\left\vert \widetilde{K}%
\right\vert =n-k$, so that $\underbrace{\left\vert K\right\vert }%
_{=k}+\underbrace{\left\vert \widetilde{K}\right\vert }_{=n-k}=k+\left(
n-k\right)  =n$. Also, $\underbrace{\left\vert \widetilde{K}\right\vert
}_{=n-k}+\underbrace{\left\vert P\right\vert }_{=k}=n-k+k=n$ and
$\underbrace{\left\vert \widetilde{Q}\right\vert }_{=n-k}%
+\underbrace{\left\vert K\right\vert }_{=k}=n-k+k=n$.
\end{vershort}

\begin{verlong}
Let $K$ be the set $\left\{  1,2,\ldots,k\right\}  $. This is a subset of
$\left\{  1,2,\ldots,n\right\}  $, since $k\in\left\{  0,1,\ldots,n\right\}
$. Its size is $\left\vert \underbrace{K}_{=\left\{  1,2,\ldots,k\right\}
}\right\vert =\left\vert \left\{  1,2,\ldots,k\right\}  \right\vert =k$. Its
complement is%
\begin{align*}
\widetilde{K}  &  =\left\{  1,2,\ldots,n\right\}  \setminus\underbrace{K}%
_{=\left\{  1,2,\ldots,k\right\}  }\ \ \ \ \ \ \ \ \ \ \left(  \text{by the
definition of }\widetilde{K}\right) \\
&  =\left\{  1,2,\ldots,n\right\}  \setminus\left\{  1,2,\ldots,k\right\}
=\left\{  k+1,k+2,\ldots,n\right\}  .
\end{align*}
Hence, $\left\vert \widetilde{K}\right\vert =\left\vert \left\{
k+1,k+2,\ldots,n\right\}  \right\vert =n-k$ (since $k\in\left\{
0,1,\ldots,n\right\}  $), so that $\underbrace{\left\vert K\right\vert }%
_{=k}+\underbrace{\left\vert \widetilde{K}\right\vert }_{=n-k}=k+\left(
n-k\right)  =n$. Also, $\underbrace{\left\vert \widetilde{K}\right\vert
}_{=n-k}+\underbrace{\left\vert P\right\vert }_{=k}=n-k+k=n$ and
$\underbrace{\left\vert \widetilde{Q}\right\vert }_{=n-k}%
+\underbrace{\left\vert K\right\vert }_{=k}=n-k+k=n$.
\end{verlong}

\begin{vershort}
Clearly, every subset $I$ of $\left\{  1,2,\ldots,n\right\}  $ satisfies
$\widetilde{\widetilde{I}}=I$. Thus, $\widetilde{\widetilde{Q}}=Q$ and
$\widetilde{\widetilde{K}}=K$.
\end{vershort}

\begin{verlong}
The definition of $\widetilde{\widetilde{Q}}$ yields%
\[
\widetilde{\widetilde{Q}}=\left\{  1,2,\ldots,n\right\}  \setminus
\underbrace{\widetilde{Q}}_{\substack{=\left\{  1,2,\ldots,n\right\}
\setminus Q\\\text{(by the definition of }\widetilde{Q}\text{)}}}=\left\{
1,2,\ldots,n\right\}  \setminus\left(  \left\{  1,2,\ldots,n\right\}
\setminus Q\right)  =Q
\]
(since $Q$ is a subset of $\left\{  1,2,\ldots,n\right\}  $). The same
argument (applied to $K$ instead of $Q$) yields $\widetilde{\widetilde{K}}=K$.
\end{verlong}

\begin{vershort}
From $K=\left\{  1,2,\ldots,k\right\}  $, we obtain $w\left(  K\right)
=\left(  1,2,\ldots,k\right)  $ (since $w\left(  K\right)  $ is defined as the
list of all elements of $K$ in increasing order (with no repetitions)).
Likewise, from $\widetilde{K}=\left\{  k+1,k+2,\ldots,n\right\}  $, we obtain
$w\left(  \widetilde{K}\right)  =\left(  k+1,k+2,\ldots,n\right)  $.
\end{vershort}

\begin{verlong}
The definition of $w\left(  K\right)  $ yields%
\begin{align*}
w\left(  K\right)   &  =\left(  \text{the list of all elements of }K\text{ in
increasing order (with no repetitions)}\right) \\
&  =\left(  \text{the list of all elements of }\left\{  1,2,\ldots,k\right\}
\text{ in increasing }\right. \\
&  \ \ \ \ \ \ \ \ \ \ \left.  \text{order (with no repetitions)}\right)
\ \ \ \ \ \ \ \ \ \ \left(  \text{since }K=\left\{  1,2,\ldots,k\right\}
\right) \\
&  =\left(  1,2,\ldots,k\right)  .
\end{align*}

The definition of $w\left(  \widetilde{K}\right)  $ yields%
\begin{align*}
w\left(  \widetilde{K}\right)   &  =\left(  \text{the list of all elements of
}\widetilde{K}\text{ in increasing order (with no repetitions)}\right) \\
&  =\left(  \text{the list of all elements of }\left\{  k+1,k+2,\ldots
,n\right\}  \text{ in increasing}\right. \\
&  \ \ \ \ \ \ \ \ \ \ \left.  \text{order (with no repetitions)}\right)
\ \ \ \ \ \ \ \ \ \ \left(  \text{since }\widetilde{K}=\left\{  k+1,k+2,\ldots
,n\right\}  \right) \\
&  =\left(  k+1,k+2,\ldots,n\right)  .
\end{align*}

\end{verlong}

We shall now prove the following nine claims:

\begin{statement}
\textit{Claim 1:} Every $i\in\widetilde{K}$ and $j\in P$ satisfy $c_{i,j}=0$.
\end{statement}

\begin{statement}
\textit{Claim 2:} We have $\operatorname*{sub}\nolimits_{w\left(  K\right)
}^{w\left(  P\right)  }C=\operatorname*{cols}\nolimits_{w\left(  P\right)  }A$.
\end{statement}

\begin{statement}
\textit{Claim 3:} We have $\operatorname*{sub}\nolimits_{w\left(
\widetilde{K}\right)  }^{w\left(  \widetilde{P}\right)  }C=I_{n-k}$.
\end{statement}

\begin{statement}
\textit{Claim 4:} Every $i\in\widetilde{Q}$ and $j\in K$ satisfy $d_{i,j}=0$.
\end{statement}

\begin{statement}
\textit{Claim 5:} We have $\operatorname*{sub}\nolimits_{w\left(
\widetilde{Q}\right)  }^{w\left(  \widetilde{K}\right)  }%
D=\operatorname*{rows}\nolimits_{w\left(  \widetilde{Q}\right)  }B$.
\end{statement}

\begin{statement}
\textit{Claim 6:} We have $\operatorname*{sub}\nolimits_{w\left(  Q\right)
}^{w\left(  K\right)  }D=I_{k}$.
\end{statement}

\begin{statement}
\textit{Claim 7:} Every $i\in K$ and $j\in\widetilde{K}$ satisfy $g_{i,j}=0$.
\end{statement}

\begin{statement}
\textit{Claim 8:} We have $\operatorname*{sub}\nolimits_{w\left(  K\right)
}^{w\left(  K\right)  }\left(  CD\right)  =\operatorname*{cols}%
\nolimits_{w\left(  Q\right)  }A$.
\end{statement}

\begin{statement}
\textit{Claim 9:} We have $\operatorname*{sub}\nolimits_{w\left(
\widetilde{K}\right)  }^{w\left(  \widetilde{K}\right)  }\left(  CD\right)
=\operatorname*{rows}\nolimits_{w\left(  \widetilde{P}\right)  }B$.
\end{statement}

Once these nine claims are proved, we will be able to factor each of the three
determinants appearing in (\ref{sol.det.AB=0-proport.detCD=}) using Exercise
\ref{exe.det.blocktria-twisted} \textbf{(b)}, and this will quickly yield the
claim of Exercise \ref{exe.det.AB=0-proport}. So let us prove the above nine
claims: \medskip

[\textit{Proof of Claim 1:} Let $i\in\widetilde{K}$ and $j\in P$. Then,
$i\in\widetilde{K}=\left\{  k+1,k+2,\ldots,n\right\}  $, so that $i\geq
k+1>k$. But the definition of $c_{i,j}$ yields%
\begin{align*}
c_{i,j}  &  =%
\begin{cases}
a_{i,j}, & \text{if }i\leq k;\\
\left[  j=\widetilde{p}_{i-k}\right]  , & \text{if }i>k
\end{cases}
\\
&  =\left[  j=\widetilde{p}_{i-k}\right]  \ \ \ \ \ \ \ \ \ \ \left(
\text{since }i>k\right)  .
\end{align*}
However, $\widetilde{p}_{i-k}$ is an entry of the list $w\left(
\widetilde{P}\right)  $ (since $w\left(  \widetilde{P}\right)  =\left(
\widetilde{p}_{1},\widetilde{p}_{2},\ldots,\widetilde{p}_{n-k}\right)  $). On
the other hand, all entries of $w\left(  \widetilde{P}\right)  $ belong to
$\widetilde{P}$ (since $w\left(  \widetilde{P}\right)  $ was defined as the
list of all elements of $\widetilde{P}$ in increasing order (with no
repetitions)). Thus, $\widetilde{p}_{i-k}$ belongs to $\widetilde{P}$ (since
$\widetilde{p}_{i-k}$ is an entry of the list $w\left(  \widetilde{P}\right)
$). Therefore,%
\[
\widetilde{p}_{i-k}\in\widetilde{P}=\left\{  1,2,\ldots,n\right\}  \setminus
P\ \ \ \ \ \ \ \ \ \ \left(  \text{by the definition of }\widetilde{P}\right)
.
\]
Thus, $\widetilde{p}_{i-k}\notin P$. Hence, we cannot have $j=\widetilde{p}%
_{i-k}$ (since this would entail $j=\widetilde{p}_{i-k}\notin P$, which would
contradict $j\in P$). Therefore, $\left[  j=\widetilde{p}_{i-k}\right]  =0$.
Hence, $c_{i,j}=\left[  j=\widetilde{p}_{i-k}\right]  =0$. This proves Claim
1.] \medskip

[\textit{Proof of Claim 2:} From $w\left(  P\right)  =\left(  p_{1}%
,p_{2},\ldots,p_{k}\right)  $ and $w\left(  K\right)  =\left(  1,2,\ldots
,k\right)  $, we obtain%
\begin{align}
\operatorname*{sub}\nolimits_{w\left(  K\right)  }^{w\left(  P\right)  }C  &
=\operatorname*{sub}\nolimits_{\left(  1,2,\ldots,k\right)  }^{\left(
p_{1},p_{2},\ldots,p_{k}\right)  }C\nonumber\\
&  =\operatorname*{sub}\nolimits_{1,2,\ldots,k}^{p_{1},p_{2},\ldots,p_{k}%
}C\ \ \ \ \ \ \ \ \ \ \left(  \text{by the definition of }\operatorname*{sub}%
\nolimits_{\left(  1,2,\ldots,k\right)  }^{\left(  p_{1},p_{2},\ldots
,p_{k}\right)  }C\right) \nonumber\\
&  =\left(  c_{x,p_{y}}\right)  _{1\leq x\leq k,\ 1\leq y\leq k}
\label{sol.det.AB=0-proport.c2.pf.1}%
\end{align}
(by Definition \ref{def.submatrix}, since $C=\left(  c_{i,j}\right)  _{1\leq
i\leq n,\ 1\leq j\leq n}$). However, for each $x\in\left\{  1,2,\ldots
,k\right\}  $ and $y\in\left\{  1,2,\ldots,k\right\}  $, we have%
\begin{align*}
c_{x,p_{y}}  &  =%
\begin{cases}
a_{x,p_{y}}, & \text{if }x\leq k;\\
\left[  p_{y}=\widetilde{p}_{x-k}\right]  , & \text{if }x>k
\end{cases}
\ \ \ \ \ \ \ \ \ \ \left(  \text{by the definition of }c_{x,p_{y}}\right) \\
&  =a_{x,p_{y}}\ \ \ \ \ \ \ \ \ \ \left(  \text{since }x\leq k\text{ (because
}x\in\left\{  1,2,\ldots,k\right\}  \text{)}\right)  .
\end{align*}
Hence, we can rewrite (\ref{sol.det.AB=0-proport.c2.pf.1}) as%
\begin{equation}
\operatorname*{sub}\nolimits_{w\left(  K\right)  }^{w\left(  P\right)
}C=\left(  a_{x,p_{y}}\right)  _{1\leq x\leq k,\ 1\leq y\leq k}.
\label{sol.det.AB=0-proport.c2.pf.2}%
\end{equation}

On the other hand, (\ref{sol.det.AB=0-proport.colswPA=}) becomes%
\[
\operatorname*{cols}\nolimits_{w\left(  P\right)  }A=\left(  a_{i,p_{y}%
}\right)  _{1\leq i\leq k,\ 1\leq y\leq k}=\left(  a_{x,p_{y}}\right)  _{1\leq
x\leq k,\ 1\leq y\leq k}%
\]
(here, we have renamed the index $\left(  i,y\right)  $ as $\left(
x,y\right)  $). Comparing this with (\ref{sol.det.AB=0-proport.c2.pf.2}), we
obtain $\operatorname*{sub}\nolimits_{w\left(  K\right)  }^{w\left(  P\right)
}C=\operatorname*{cols}\nolimits_{w\left(  P\right)  }A$. This proves Claim
2.] \medskip

[\textit{Proof of Claim 3:} Recall that $w\left(  \widetilde{P}\right)  $ was
defined as the list of all elements of $\widetilde{P}$ in increasing order
(with no repetitions). Hence, the list $w\left(  \widetilde{P}\right)  $ has
no repetitions, i.e., no two equal entries. In other words, the list $\left(
\widetilde{p}_{1},\widetilde{p}_{2},\ldots,\widetilde{p}_{n-k}\right)  $ has
no two equal entries (since $w\left(  \widetilde{P}\right)  =\left(
\widetilde{p}_{1},\widetilde{p}_{2},\ldots,\widetilde{p}_{n-k}\right)  $). In
other words, two elements $y$ and $x$ of $\left\{  1,2,\ldots,n-k\right\}  $
satisfy $\widetilde{p}_{y}=\widetilde{p}_{x}$ if and only if they satisfy
$y=x$. In other words, for any two elements $y$ and $x$ of $\left\{
1,2,\ldots,n-k\right\}  $, the two logical statements \textquotedblleft%
$\widetilde{p}_{y}=\widetilde{p}_{x}$\textquotedblright\ and \textquotedblleft%
$y=x$\textquotedblright\ are equivalent. Hence, for any two elements $y$ and
$x$ of $\left\{  1,2,\ldots,n-k\right\}  $, we have%
\begin{equation}
\left[  \widetilde{p}_{y}=\widetilde{p}_{x}\right]  =\left[  y=x\right]
\label{sol.det.AB=0-proport.c3.pf.0}%
\end{equation}
(by Exercise \ref{exe.iverson-prop} \textbf{(a)}, applied to $\mathcal{A}%
=\left(  \text{\textquotedblleft}\widetilde{p}_{y}=\widetilde{p}%
_{x}\text{\textquotedblright}\right)  $ and $\mathcal{B}=\left(
\text{\textquotedblleft}y=x\text{\textquotedblright}\right)  $).

From $w\left(  \widetilde{P}\right)  =\left(  \widetilde{p}_{1},\widetilde{p}%
_{2},\ldots,\widetilde{p}_{n-k}\right)  $ and $w\left(  \widetilde{K}\right)
=\left(  k+1,k+2,\ldots,n\right)  $, we obtain%
\begin{align}
\operatorname*{sub}\nolimits_{w\left(  \widetilde{K}\right)  }^{w\left(
\widetilde{P}\right)  }C  &  =\operatorname*{sub}\nolimits_{\left(
k+1,k+2,\ldots,n\right)  }^{\left(  \widetilde{p}_{1},\widetilde{p}_{2}%
,\ldots,\widetilde{p}_{n-k}\right)  }C\nonumber\\
&  =\operatorname*{sub}\nolimits_{k+1,k+2,\ldots,n}^{\widetilde{p}%
_{1},\widetilde{p}_{2},\ldots,\widetilde{p}_{n-k}}C\ \ \ \ \ \ \ \ \ \ \left(
\text{by the definition of }\operatorname*{sub}\nolimits_{\left(
k+1,k+2,\ldots,n\right)  }^{\left(  \widetilde{p}_{1},\widetilde{p}_{2}%
,\ldots,\widetilde{p}_{n-k}\right)  }C\right) \nonumber\\
&  =\left(  c_{k+x,\widetilde{p}_{y}}\right)  _{1\leq x\leq n-k,\ 1\leq y\leq
n-k} \label{sol.det.AB=0-proport.c3.pf.1}%
\end{align}
(by Definition \ref{def.submatrix}, since $C=\left(  c_{i,j}\right)  _{1\leq
i\leq n,\ 1\leq j\leq n}$). However, for each $x\in\left\{  1,2,\ldots
,n-k\right\}  $ and $y\in\left\{  1,2,\ldots,n-k\right\}  $, we have%
\begin{align*}
c_{k+x,\widetilde{p}_{y}}  &  =%
\begin{cases}
a_{k+x,\widetilde{p}_{y}}, & \text{if }k+x\leq k;\\
\left[  \widetilde{p}_{y}=\widetilde{p}_{\left(  k+x\right)  -k}\right]  , &
\text{if }k+x>k
\end{cases}
\ \ \ \ \ \ \ \ \ \ \left(  \text{by the definition of }c_{k+x,\widetilde{p}%
_{y}}\right) \\
&  =\left[  \widetilde{p}_{y}=\widetilde{p}_{\left(  k+x\right)  -k}\right]
\ \ \ \ \ \ \ \ \ \ \left(  \text{since }k+\underbrace{x}_{\geq1>0}>k\right)
\\
&  =\left[  \widetilde{p}_{y}=\widetilde{p}_{x}\right]
\ \ \ \ \ \ \ \ \ \ \left(  \text{since }\left(  k+x\right)  -k=x\right) \\
&  =\left[  y=x\right]  \ \ \ \ \ \ \ \ \ \ \left(  \text{by
(\ref{sol.det.AB=0-proport.c3.pf.0})}\right) \\
&  =\left[  x=y\right]  \ \ \ \ \ \ \ \ \ \ \left(  \text{since the statement
\textquotedblleft}y=x\text{\textquotedblright\ is clearly equivalent to
\textquotedblleft}x=y\text{\textquotedblright}\right)  .
\end{align*}
Hence, we can rewrite (\ref{sol.det.AB=0-proport.c3.pf.1}) as%
\[
\operatorname*{sub}\nolimits_{w\left(  \widetilde{K}\right)  }^{w\left(
\widetilde{P}\right)  }C=\left(  \left[  x=y\right]  \right)  _{1\leq x\leq
n-k,\ 1\leq y\leq n-k}=\left(  \left[  i=j\right]  \right)  _{1\leq i\leq
n-k,\ 1\leq j\leq n-k}%
\]
(here, we have renamed the index $\left(  x,y\right)  $ as $\left(
i,j\right)  $). Comparing this with (\ref{sol.det.AB=0-proport.In-k=}), we
obtain $\operatorname*{sub}\nolimits_{w\left(  \widetilde{K}\right)
}^{w\left(  \widetilde{P}\right)  }C=I_{n-k}$. This proves Claim 3.]\medskip

[\textit{Proof of Claim 4:} Let $i\in\widetilde{Q}$ and $j\in K$. Then, $j\in
K=\left\{  1,2,\ldots,k\right\}  $, so that $j\leq k$. But the definition of
$d_{i,j}$ yields%
\begin{align*}
d_{i,j}  &  =%
\begin{cases}
\left[  i=q_{j}\right]  , & \text{if }j\leq k;\\
b_{i,j-k}, & \text{if }j>k
\end{cases}
\\
&  =\left[  i=q_{j}\right]  \ \ \ \ \ \ \ \ \ \ \left(  \text{since }j\leq
k\right)  .
\end{align*}
However, $q_{j}$ is an entry of the list $w\left(  Q\right)  $ (since
$w\left(  Q\right)  =\left(  q_{1},q_{2},\ldots,q_{k}\right)  $). On the other
hand, all entries of $w\left(  Q\right)  $ belong to $Q$ (since $w\left(
Q\right)  $ was defined as the list of all elements of $Q$ in increasing order
(with no repetitions)). Thus, $q_{j}$ belongs to $Q$ (since $q_{j}$ is an
entry of the list $w\left(  Q\right)  $). In other words, $q_{j}\in Q$. On the
other hand,
\[
i\in\widetilde{Q}=\left\{  1,2,\ldots,n\right\}  \setminus
Q\ \ \ \ \ \ \ \ \ \ \left(  \text{by the definition of }\widetilde{Q}\right)
.
\]
Thus, $i\notin Q$. Hence, we cannot have $i=q_{j}$ (since this would entail
$i=q_{j}\in Q$, which would contradict $i\notin Q$). Therefore, $\left[
i=q_{j}\right]  =0$. Hence, $d_{i,j}=\left[  i=q_{j}\right]  =0$. This proves
Claim 4.] \medskip

[\textit{Proof of Claim 5:} From $w\left(  \widetilde{Q}\right)  =\left(
\widetilde{q}_{1},\widetilde{q}_{2},\ldots,\widetilde{q}_{n-k}\right)  $ and
$w\left(  \widetilde{K}\right)  =\left(  k+1,k+2,\ldots,n\right)  $, we obtain%
\begin{align}
\operatorname*{sub}\nolimits_{w\left(  \widetilde{Q}\right)  }^{w\left(
\widetilde{K}\right)  }D  &  =\operatorname*{sub}\nolimits_{\left(
\widetilde{q}_{1},\widetilde{q}_{2},\ldots,\widetilde{q}_{n-k}\right)
}^{\left(  k+1,k+2,\ldots,n\right)  }D\nonumber\\
&  =\operatorname*{sub}\nolimits_{\widetilde{q}_{1},\widetilde{q}_{2}%
,\ldots,\widetilde{q}_{n-k}}^{k+1,k+2,\ldots,n}D\ \ \ \ \ \ \ \ \ \ \left(
\text{by the definition of }\operatorname*{sub}\nolimits_{\left(
\widetilde{q}_{1},\widetilde{q}_{2},\ldots,\widetilde{q}_{n-k}\right)
}^{\left(  k+1,k+2,\ldots,n\right)  }D\right) \nonumber\\
&  =\left(  d_{\widetilde{q}_{x},k+y}\right)  _{1\leq x\leq n-k,\ 1\leq y\leq
n-k} \label{sol.det.AB=0-proport.c5.pf.1}%
\end{align}
(by Definition \ref{def.submatrix}, since $D=\left(  d_{i,j}\right)  _{1\leq
i\leq n,\ 1\leq j\leq n}$). However, for each $x\in\left\{  1,2,\ldots
,n-k\right\}  $ and $y\in\left\{  1,2,\ldots,n-k\right\}  $, we have%
\begin{align*}
d_{\widetilde{q}_{x},k+y}  &  =%
\begin{cases}
\left[  \widetilde{q}_{x}=q_{k+y}\right]  , & \text{if }k+y\leq k;\\
b_{\widetilde{q}_{x},\left(  k+y\right)  -k}, & \text{if }k+y>k
\end{cases}
\ \ \ \ \ \ \ \ \ \ \left(  \text{by the definition of }d_{\widetilde{q}%
_{x},k+y}\right) \\
&  =b_{\widetilde{q}_{x},\left(  k+y\right)  -k}\ \ \ \ \ \ \ \ \ \ \left(
\text{since }k+\underbrace{y}_{\geq1>0}>k\right) \\
&  =b_{\widetilde{q}_{x},y}\ \ \ \ \ \ \ \ \ \ \left(  \text{since }\left(
k+y\right)  -k=y\right)  .
\end{align*}
Hence, we can rewrite (\ref{sol.det.AB=0-proport.c5.pf.1}) as%
\[
\operatorname*{sub}\nolimits_{w\left(  \widetilde{Q}\right)  }^{w\left(
\widetilde{K}\right)  }D=\left(  b_{\widetilde{q}_{x},y}\right)  _{1\leq x\leq
n-k,\ 1\leq y\leq n-k}.
\]

\begin{vershort}
\noindent On the other hand, (\ref{sol.det.AB=0-proport.rowswQB=.short})
becomems%
\[
\operatorname*{rows}\nolimits_{w\left(  \widetilde{Q}\right)  }B=\left(
b_{\widetilde{q}_{x},j}\right)  _{1\leq x\leq n-k,\ 1\leq j\leq n-k}=\left(
b_{\widetilde{q}_{x},y}\right)  _{1\leq x\leq n-k,\ 1\leq y\leq n-k}%
\]
(here, we have renamed the index $\left(  x,j\right)  $ as $\left(
x,y\right)  $).
\end{vershort}

\begin{verlong}
\noindent On the other hand,%
\begin{align*}
\operatorname*{rows}\nolimits_{w\left(  \widetilde{Q}\right)  }B  &  =\left(
b_{\widetilde{q}_{x},j}\right)  _{1\leq x\leq n-k,\ 1\leq j\leq n-k}%
\ \ \ \ \ \ \ \ \ \ \left(  \text{by (\ref{sol.det.AB=0-proport.rowswQB=}%
)}\right) \\
&  =\left(  b_{\widetilde{q}_{x},y}\right)  _{1\leq x\leq n-k,\ 1\leq y\leq
n-k}\\
&  \ \ \ \ \ \ \ \ \ \ \ \ \ \ \ \ \ \ \ \ \left(  \text{here, we have renamed
the index }\left(  x,j\right)  \text{ as }\left(  x,y\right)  \right)  ,
\end{align*}

\end{verlong}

\noindent Comparing these two equalities, we obtain $\operatorname*{sub}%
\nolimits_{w\left(  \widetilde{Q}\right)  }^{w\left(  \widetilde{K}\right)
}D=\operatorname*{rows}\nolimits_{w\left(  \widetilde{Q}\right)  }B$. This
proves Claim 5.] \medskip

[\textit{Proof of Claim 6:} Recall that $w\left(  Q\right)  $ was defined as
the list of all elements of $Q$ in increasing order (with no repetitions).
Hence, the list $w\left(  Q\right)  $ has no repetitions, i.e., no two equal
entries. In other words, the list $\left(  q_{1},q_{2},\ldots,q_{k}\right)  $
has no two equal entries (since $w\left(  Q\right)  =\left(  q_{1}%
,q_{2},\ldots,q_{k}\right)  $). In other words, two elements $x$ and $y$ of
$\left\{  1,2,\ldots,k\right\}  $ satisfy $q_{x}=q_{y}$ if and only if they
satisfy $x=y$. In other words, for any two elements $y$ and $x$ of $\left\{
1,2,\ldots,k\right\}  $, the two logical statements \textquotedblleft%
$q_{x}=q_{y}$\textquotedblright\ and \textquotedblleft$x=y$\textquotedblright%
\ are equivalent. Hence, for any two elements $x$ and $y$ of $\left\{
1,2,\ldots,n-k\right\}  $, we have%
\begin{equation}
\left[  q_{x}=q_{y}\right]  =\left[  x=y\right]
\label{sol.det.AB=0-proport.c6.pf.0}%
\end{equation}
(by Exercise \ref{exe.iverson-prop} \textbf{(a)}, applied to $\mathcal{A}%
=\left(  \text{\textquotedblleft}q_{x}=q_{y}\text{\textquotedblright}\right)
$ and $\mathcal{B}=\left(  \text{\textquotedblleft}x=y\text{\textquotedblright%
}\right)  $).

From $w\left(  Q\right)  =\left(  q_{1},q_{2},\ldots,q_{k}\right)  $ and
$w\left(  K\right)  =\left(  1,2,\ldots,k\right)  $, we obtain%
\begin{align}
\operatorname*{sub}\nolimits_{w\left(  Q\right)  }^{w\left(  K\right)  }D  &
=\operatorname*{sub}\nolimits_{\left(  q_{1},q_{2},\ldots,q_{k}\right)
}^{\left(  1,2,\ldots,k\right)  }D\nonumber\\
&  =\operatorname*{sub}\nolimits_{q_{1},q_{2},\ldots,q_{k}}^{1,2,\ldots
,k}D\ \ \ \ \ \ \ \ \ \ \left(  \text{by the definition of }%
\operatorname*{sub}\nolimits_{\left(  q_{1},q_{2},\ldots,q_{k}\right)
}^{\left(  1,2,\ldots,k\right)  }D\right) \nonumber\\
&  =\left(  d_{q_{x},y}\right)  _{1\leq x\leq k,\ 1\leq y\leq k}
\label{sol.det.AB=0-proport.c6.pf.1}%
\end{align}
(by Definition \ref{def.submatrix}, since $D=\left(  d_{i,j}\right)  _{1\leq
i\leq n,\ 1\leq j\leq n}$). However, for each $x\in\left\{  1,2,\ldots
,k\right\}  $ and $y\in\left\{  1,2,\ldots,k\right\}  $, we have%
\begin{align*}
d_{q_{x},y}  &  =%
\begin{cases}
\left[  q_{x}=q_{y}\right]  , & \text{if }y\leq k;\\
b_{q_{x},y-k}, & \text{if }y>k
\end{cases}
\ \ \ \ \ \ \ \ \ \ \left(  \text{by the definition of }d_{q_{x},y}\right) \\
&  =\left[  q_{x}=q_{y}\right]  \ \ \ \ \ \ \ \ \ \ \left(  \text{since }y\leq
k\text{ (because }y\in\left\{  1,2,\ldots,k\right\}  \text{)}\right) \\
&  =\left[  x=y\right]  \ \ \ \ \ \ \ \ \ \ \left(  \text{by
(\ref{sol.det.AB=0-proport.c6.pf.0})}\right)  .
\end{align*}
Hence, we can rewrite (\ref{sol.det.AB=0-proport.c6.pf.1}) as%
\[
\operatorname*{sub}\nolimits_{w\left(  Q\right)  }^{w\left(  K\right)
}D=\left(  \left[  x=y\right]  \right)  _{1\leq x\leq k,\ 1\leq y\leq
k}=\left(  \left[  i=j\right]  \right)  _{1\leq i\leq k,\ 1\leq j\leq k}%
\]
(here, we have renamed the index $\left(  x,y\right)  $ as $\left(
i,j\right)  $). Comparing this with (\ref{sol.det.AB=0-proport.Ik=}), we
obtain $\operatorname*{sub}\nolimits_{w\left(  Q\right)  }^{w\left(  K\right)
}D=I_{k}$. This proves Claim 6.]\medskip

[\textit{Proof of Claim 7:} Let $i\in K$ and $j\in\widetilde{K}$. Then,
$j\in\widetilde{K}=\left\{  k+1,k+2,\ldots,n\right\}  $, so that $j\geq
k+1>k$. Also, $i\in K=\left\{  1,2,\ldots,k\right\}  $, so that $i\leq k$.
From $j\in\left\{  k+1,k+2,\ldots,n\right\}  $, we obtain $j-k\in\left\{
1,2,\ldots,n-k\right\}  $.

Now, for each $m\in\left\{  1,2,\ldots,n\right\}  $, we have%
\begin{align}
c_{i,m}  &  =%
\begin{cases}
a_{i,m}, & \text{if }i\leq k;\\
\left[  m=\widetilde{p}_{i-k}\right]  , & \text{if }i>k
\end{cases}
\ \ \ \ \ \ \ \ \ \ \left(  \text{by the definition of }c_{i,m}\right)
\nonumber\\
&  =a_{i,m}\ \ \ \ \ \ \ \ \ \ \left(  \text{since }i\leq k\right)
\label{sol.det.AB=0-proport.c7.pf.1}%
\end{align}
and%
\begin{align}
d_{m,j}  &  =%
\begin{cases}
\left[  m=q_{j}\right]  , & \text{if }j\leq k;\\
b_{m,j-k}, & \text{if }j>k
\end{cases}
\ \ \ \ \ \ \ \ \ \ \left(  \text{by the definition of }d_{m,j}\right)
\nonumber\\
&  =b_{m,j-k}\ \ \ \ \ \ \ \ \ \ \left(  \text{since }j>k\right)  .
\label{sol.det.AB=0-proport.c7.pf.2}%
\end{align}
Now, (\ref{sol.det.AB=0-proport.gij=}) yields%
\[
g_{i,j}=\sum_{m=1}^{n}\underbrace{c_{i,m}}_{\substack{=a_{i,m}\\\text{(by
(\ref{sol.det.AB=0-proport.c7.pf.1}))}}}\ \ \underbrace{d_{m,j}}%
_{\substack{=b_{m,j-k}\\\text{(by (\ref{sol.det.AB=0-proport.c7.pf.2}))}%
}}=\sum_{m=1}^{n}a_{i,m}b_{m,j-k}.
\]
On the other hand, (\ref{sol.det.AB=0-proport.0=sumab}) (applied to $j-k$
instead of $j$) yields%
\[
0=\sum_{m=1}^{n}a_{i,m}b_{m,j-k}.
\]
Comparing these two equalities, we obtain $g_{i,j}=0$. This proves Claim 7.]
\medskip

[\textit{Proof of Claim 8:} From $w\left(  K\right)  =\left(  1,2,\ldots
,k\right)  $, we obtain%
\begin{align*}
\operatorname*{sub}\nolimits_{w\left(  K\right)  }^{w\left(  K\right)
}\left(  CD\right)   &  =\operatorname*{sub}\nolimits_{\left(  1,2,\ldots
,k\right)  }^{\left(  1,2,\ldots,k\right)  }\left(  CD\right) \\
&  =\operatorname*{sub}\nolimits_{1,2,\ldots,k}^{1,2,\ldots,k}\left(
CD\right)  \ \ \ \ \ \ \ \ \ \ \left(  \text{by the definition of
}\operatorname*{sub}\nolimits_{\left(  1,2,\ldots,k\right)  }^{\left(
1,2,\ldots,k\right)  }\left(  CD\right)  \right) \\
&  =\left(  g_{x,y}\right)  _{1\leq x\leq k,\ 1\leq y\leq k}%
\end{align*}
(by Definition \ref{def.submatrix}, since $CD=\left(  g_{i,j}\right)  _{1\leq
i\leq n,\ 1\leq j\leq n}$). Thus,%
\begin{align}
\operatorname*{sub}\nolimits_{w\left(  K\right)  }^{w\left(  K\right)
}\left(  CD\right)   &  =\left(  g_{x,y}\right)  _{1\leq x\leq k,\ 1\leq y\leq
k}\nonumber\\
&  =\left(  g_{i,j}\right)  _{1\leq i\leq k,\ 1\leq j\leq k}
\label{sol.det.AB=0-proport.c8.pf.1}%
\end{align}
(here, we have renamed the index $\left(  x,y\right)  $ as $\left(
i,j\right)  $).

Now, let $i\in\left\{  1,2,\ldots,k\right\}  $ and $j\in\left\{
1,2,\ldots,k\right\}  $ be arbitrary. Thus, $i\leq k$ and $j\leq k$. For each
$m\in\left\{  1,2,\ldots,n\right\}  $, we have%
\begin{align}
c_{i,m}  &  =%
\begin{cases}
a_{i,m}, & \text{if }i\leq k;\\
\left[  m=\widetilde{p}_{i-k}\right]  , & \text{if }i>k
\end{cases}
\ \ \ \ \ \ \ \ \ \ \left(  \text{by the definition of }c_{i,m}\right)
\nonumber\\
&  =a_{i,m}\ \ \ \ \ \ \ \ \ \ \left(  \text{since }i\leq k\right)
\label{sol.det.AB=0-proport.c8.pf.2}%
\end{align}
and%
\begin{align}
d_{m,j}  &  =%
\begin{cases}
\left[  m=q_{j}\right]  , & \text{if }j\leq k;\\
b_{m,j-k}, & \text{if }j>k
\end{cases}
\ \ \ \ \ \ \ \ \ \ \left(  \text{by the definition of }d_{m,j}\right)
\nonumber\\
&  =\left[  m=q_{j}\right]  \ \ \ \ \ \ \ \ \ \ \left(  \text{since }j\leq
k\right)  . \label{sol.det.AB=0-proport.c8.pf.3}%
\end{align}
Now, (\ref{sol.det.AB=0-proport.gij=}) yields%
\begin{align*}
g_{i,j}  &  =\underbrace{\sum_{m=1}^{n}}_{=\sum_{m\in\left\{  1,2,\ldots
,n\right\}  }}\underbrace{c_{i,m}}_{\substack{=a_{i,m}\\\text{(by
(\ref{sol.det.AB=0-proport.c8.pf.2}))}}}\ \ \underbrace{d_{m,j}}%
_{\substack{=\left[  m=q_{j}\right]  \\\text{(by
(\ref{sol.det.AB=0-proport.c8.pf.3}))}}}\\
&  =\sum_{m\in\left\{  1,2,\ldots,n\right\}  }a_{i,m}\left[  m=q_{j}\right]
=a_{i,q_{j}}\underbrace{\left[  q_{j}=q_{j}\right]  }%
_{\substack{=1\\\text{(since }q_{j}=q_{j}\text{)}}}+\sum_{\substack{m\in
\left\{  1,2,\ldots,n\right\}  ;\\m\neq q_{j}}}a_{i,m}\underbrace{\left[
m=q_{j}\right]  }_{\substack{=0\\\text{(since }m\neq q_{j}\text{)}}}\\
&  \ \ \ \ \ \ \ \ \ \ \ \ \ \ \ \ \ \ \ \ \left(  \text{here, we have split
off the addend for }m=q_{j}\text{ from the sum}\right) \\
&  =a_{i,q_{j}}+\underbrace{\sum_{\substack{m\in\left\{  1,2,\ldots,n\right\}
;\\m\neq q_{j}}}a_{i,m}0}_{=0}=a_{i,q_{j}}.
\end{align*}

Forget that we fixed $i$ and $j$. We thus have proved that%
\[
g_{i,j}=a_{i,q_{j}}\ \ \ \ \ \ \ \ \ \ \text{for all }i\in\left\{
1,2,\ldots,k\right\}  \text{ and }j\in\left\{  1,2,\ldots,k\right\}  .
\]
Hence, we can rewrite (\ref{sol.det.AB=0-proport.c8.pf.1}) as follows:%
\begin{equation}
\operatorname*{sub}\nolimits_{w\left(  K\right)  }^{w\left(  K\right)
}\left(  CD\right)  =\left(  a_{i,q_{j}}\right)  _{1\leq i\leq k,\ 1\leq j\leq
k}=\left(  a_{i,q_{y}}\right)  _{1\leq i\leq k,\ 1\leq y\leq k}\nonumber
\end{equation}
(here, we have renamed the index $\left(  i,j\right)  $ as $\left(
i,y\right)  $).

\begin{vershort}
Comparing this with (\ref{sol.det.AB=0-proport.colswQA=.short}), we obtain
$\operatorname*{sub}\nolimits_{w\left(  K\right)  }^{w\left(  K\right)
}\left(  CD\right)  =\operatorname*{cols}\nolimits_{w\left(  Q\right)  }A$.
\end{vershort}

\begin{verlong}
Comparing this with (\ref{sol.det.AB=0-proport.colswQA=}), we obtain
$\operatorname*{sub}\nolimits_{w\left(  K\right)  }^{w\left(  K\right)
}\left(  CD\right)  =\operatorname*{cols}\nolimits_{w\left(  Q\right)  }A$.
\end{verlong}

This proves Claim 8.] \medskip

[\textit{Proof of Claim 9:} From $w\left(  \widetilde{K}\right)  =\left(
k+1,k+2,\ldots,n\right)  $, we obtain%
\begin{align}
\operatorname*{sub}\nolimits_{w\left(  \widetilde{K}\right)  }^{w\left(
\widetilde{K}\right)  }\left(  CD\right)   &  =\operatorname*{sub}%
\nolimits_{\left(  k+1,k+2,\ldots,n\right)  }^{\left(  k+1,k+2,\ldots
,n\right)  }\left(  CD\right) \nonumber\\
&  =\operatorname*{sub}\nolimits_{k+1,k+2,\ldots,n}^{k+1,k+2,\ldots,n}\left(
CD\right)  \ \ \ \ \ \ \ \ \ \ \left(  \text{by the definition of
}\operatorname*{sub}\nolimits_{\left(  k+1,k+2,\ldots,n\right)  }^{\left(
k+1,k+2,\ldots,n\right)  }\left(  CD\right)  \right) \nonumber\\
&  =\left(  g_{k+x,k+y}\right)  _{1\leq x\leq n-k,\ 1\leq y\leq n-k}
\label{sol.det.AB=0-proport.c9.pf.0}%
\end{align}
(by Definition \ref{def.submatrix}, since $CD=\left(  g_{i,j}\right)  _{1\leq
i\leq n,\ 1\leq j\leq n}$).

Now, let $x\in\left\{  1,2,\ldots,n-k\right\}  $ and $y\in\left\{
1,2,\ldots,n-k\right\}  $ be arbitrary. Set $i=k+x$ and $j=k+y$. Then,
$i=k+x\in\left\{  k+1,k+2,\ldots,n\right\}  $ (since $x\in\left\{
1,2,\ldots,n-k\right\}  $), so that $i\geq k+1>k$. Also, $j=k+y\in\left\{
k+1,k+2,\ldots,n\right\}  $ (since $y\in\left\{  1,2,\ldots,n-k\right\}  $),
so that $j\geq k+1>k$. Furthermore, $i-k=x$ (since $i=k+x$) and $j-k=y$ (since
$j=k+y$). For each $m\in\left\{  1,2,\ldots,n\right\}  $, we have%
\begin{align}
c_{i,m}  &  =%
\begin{cases}
a_{i,m}, & \text{if }i\leq k;\\
\left[  m=\widetilde{p}_{i-k}\right]  , & \text{if }i>k
\end{cases}
\ \ \ \ \ \ \ \ \ \ \left(  \text{by the definition of }c_{i,m}\right)
\nonumber\\
&  =\left[  m=\widetilde{p}_{i-k}\right]  \ \ \ \ \ \ \ \ \ \ \left(
\text{since }i>k\right) \nonumber\\
&  =\left[  m=\widetilde{p}_{x}\right]  \ \ \ \ \ \ \ \ \ \ \left(
\text{since }i-k=x\right)  \label{sol.det.AB=0-proport.c9.pf.1}%
\end{align}
and%
\begin{align}
d_{m,j}  &  =%
\begin{cases}
\left[  m=q_{j}\right]  , & \text{if }j\leq k;\\
b_{m,j-k}, & \text{if }j>k
\end{cases}
\ \ \ \ \ \ \ \ \ \ \left(  \text{by the definition of }d_{m,j}\right)
\nonumber\\
&  =b_{m,j-k}\ \ \ \ \ \ \ \ \ \ \left(  \text{since }j>k\right) \nonumber\\
&  =b_{m,y}\ \ \ \ \ \ \ \ \ \ \left(  \text{since }j-k=y\right)  .
\label{sol.det.AB=0-proport.c9.pf.2}%
\end{align}
Now, (\ref{sol.det.AB=0-proport.gij=}) yields%
\begin{align*}
g_{i,j}  &  =\underbrace{\sum_{m=1}^{n}}_{=\sum_{m\in\left\{  1,2,\ldots
,n\right\}  }}\underbrace{c_{i,m}}_{\substack{=\left[  m=\widetilde{p}%
_{x}\right]  \\\text{(by (\ref{sol.det.AB=0-proport.c9.pf.1}))}}%
}\ \ \underbrace{d_{m,j}}_{\substack{=b_{m,y}\\\text{(by
(\ref{sol.det.AB=0-proport.c9.pf.2}))}}}\\
&  =\sum_{m\in\left\{  1,2,\ldots,n\right\}  }\left[  m=\widetilde{p}%
_{x}\right]  b_{m,y}=\underbrace{\left[  \widetilde{p}_{x}=\widetilde{p}%
_{x}\right]  }_{\substack{=1\\\text{(since }\widetilde{p}_{x}=\widetilde{p}%
_{x}\text{)}}}b_{\widetilde{p}_{x},y}+\sum_{\substack{m\in\left\{
1,2,\ldots,n\right\}  ;\\m\neq\widetilde{p}_{x}}}\underbrace{\left[
m=\widetilde{p}_{x}\right]  }_{\substack{=0\\\text{(since }m\neq
\widetilde{p}_{x}\text{)}}}b_{m,y}\\
&  \ \ \ \ \ \ \ \ \ \ \ \ \ \ \ \ \ \ \ \ \left(  \text{here, we have split
off the addend for }m=\widetilde{p}_{x}\text{ from the sum}\right) \\
&  =b_{\widetilde{p}_{x},y}+\underbrace{\sum_{\substack{m\in\left\{
1,2,\ldots,n\right\}  ;\\m\neq\widetilde{p}_{x}}}0b_{m,y}}_{=0}%
=b_{\widetilde{p}_{x},y}.
\end{align*}
In view of $i=k+x$ and $j=k+y$, we can rewrite this as $g_{k+x,k+y}%
=b_{\widetilde{p}_{x},y}$.

Forget that we fixed $x$ and $y$. We thus have proved that%
\[
g_{k+x,k+y}=b_{\widetilde{p}_{x},y}\ \ \ \ \ \ \ \ \ \ \text{for all }%
x\in\left\{  1,2,\ldots,n-k\right\}  \text{ and }y\in\left\{  1,2,\ldots
,n-k\right\}  .
\]
Hence, we can rewrite (\ref{sol.det.AB=0-proport.c9.pf.0}) as follows:%
\begin{equation}
\operatorname*{sub}\nolimits_{w\left(  \widetilde{K}\right)  }^{w\left(
\widetilde{K}\right)  }\left(  CD\right)  =\left(  b_{\widetilde{p}_{x}%
,y}\right)  _{1\leq x\leq n-k,\ 1\leq y\leq n-k}.\nonumber
\end{equation}
Comparing this with%
\begin{align*}
\operatorname*{rows}\nolimits_{w\left(  \widetilde{P}\right)  }B  &  =\left(
b_{\widetilde{p}_{x},j}\right)  _{1\leq x\leq n-k,\ 1\leq j\leq n-k}%
\ \ \ \ \ \ \ \ \ \ \left(  \text{by (\ref{sol.det.AB=0-proport.rowswPB=}%
)}\right) \\
&  =\left(  b_{\widetilde{p}_{x},y}\right)  _{1\leq x\leq n-k,\ 1\leq y\leq
n-k}\ \ \ \ \ \ \ \ \ \ \left(
\begin{array}
[c]{c}%
\text{here, we have renamed the}\\
\text{index }\left(  x,j\right)  \text{ as }\left(  x,y\right)
\end{array}
\right)  ,
\end{align*}
we obtain $\operatorname*{sub}\nolimits_{w\left(  \widetilde{K}\right)
}^{w\left(  \widetilde{K}\right)  }\left(  CD\right)  =\operatorname*{rows}%
\nolimits_{w\left(  \widetilde{P}\right)  }B$. This proves Claim 9.] \medskip

Having proved all nine claims, we are now ready for the endgame.

We know that $\widetilde{K}$ and $P$ are two subsets of $\left\{
1,2,\ldots,n\right\}  $ satisfying $\left\vert \widetilde{K}\right\vert
+\left\vert P\right\vert =n$. Also, $C=\left(  c_{i,j}\right)  _{1\leq i\leq
n,\ 1\leq j\leq n}\in\mathbb{K}^{n\times n}$ is an $n\times n$-matrix such
that%
\[
\text{every }i\in\widetilde{K}\text{ and }j\in P\text{ satisfy }c_{i,j}=0
\]
(by Claim 1). Hence, Exercise \ref{exe.det.blocktria-twisted} \textbf{(b)}
(applied to $\widetilde{K}$, $P$, $C$ and $c_{i,j}$ instead of $P$, $Q$, $A$
and $a_{i,j}$) yields
\begin{align}
\det C  &  =\left(  -1\right)  ^{\sum\widetilde{K}+\sum\widetilde{P}}%
\det\left(  \operatorname*{sub}\nolimits_{w\left(  \widetilde{K}\right)
}^{w\left(  \widetilde{P}\right)  }C\right)  \det\left(  \operatorname*{sub}%
\nolimits_{w\left(  \widetilde{\widetilde{K}}\right)  }^{w\left(  P\right)
}C\right) \nonumber\\
&  =\left(  -1\right)  ^{\sum\widetilde{K}+\sum\widetilde{P}}\det
\underbrace{\left(  \operatorname*{sub}\nolimits_{w\left(  \widetilde{K}%
\right)  }^{w\left(  \widetilde{P}\right)  }C\right)  }_{\substack{=I_{n-k}%
\\\text{(by Claim 3)}}}\det\underbrace{\left(  \operatorname*{sub}%
\nolimits_{w\left(  K\right)  }^{w\left(  P\right)  }C\right)  }%
_{\substack{=\operatorname*{cols}\nolimits_{w\left(  P\right)  }A\\\text{(by
Claim 2)}}}\ \ \ \ \ \ \ \ \ \ \left(  \text{since }\widetilde{\widetilde{K}%
}=K\right) \nonumber\\
&  =\left(  -1\right)  ^{\sum\widetilde{K}+\sum\widetilde{P}}\underbrace{\det
\left(  I_{n-k}\right)  }_{=1}\det\left(  \operatorname*{cols}%
\nolimits_{w\left(  P\right)  }A\right) \nonumber\\
&  =\left(  -1\right)  ^{\sum\widetilde{K}+\sum\widetilde{P}}\det\left(
\operatorname*{cols}\nolimits_{w\left(  P\right)  }A\right)  .
\label{sol.det.AB=0-proport.detC=}%
\end{align}

We know that $\widetilde{Q}$ and $K$ are two subsets of $\left\{
1,2,\ldots,n\right\}  $ satisfying $\left\vert \widetilde{Q}\right\vert
+\left\vert K\right\vert =n$. Also, $D=\left(  d_{i,j}\right)  _{1\leq i\leq
n,\ 1\leq j\leq n}\in\mathbb{K}^{n\times n}$ is an $n\times n$-matrix such
that%
\[
\text{every }i\in\widetilde{Q}\text{ and }j\in K\text{ satisfy }d_{i,j}=0
\]
(by Claim 4). Hence, Exercise \ref{exe.det.blocktria-twisted} \textbf{(b)}
(applied to $\widetilde{Q}$, $K$, $D$ and $d_{i,j}$ instead of $P$, $Q$, $A$
and $a_{i,j}$) yields%
\begin{align}
\det D  &  =\left(  -1\right)  ^{\sum\widetilde{Q}+\sum\widetilde{K}}%
\det\left(  \operatorname*{sub}\nolimits_{w\left(  \widetilde{Q}\right)
}^{w\left(  \widetilde{K}\right)  }D\right)  \det\left(  \operatorname*{sub}%
\nolimits_{w\left(  \widetilde{\widetilde{Q}}\right)  }^{w\left(  K\right)
}D\right) \nonumber\\
&  =\left(  -1\right)  ^{\sum\widetilde{Q}+\sum\widetilde{K}}\det
\underbrace{\left(  \operatorname*{sub}\nolimits_{w\left(  \widetilde{Q}%
\right)  }^{w\left(  \widetilde{K}\right)  }D\right)  }%
_{\substack{=\operatorname*{rows}\nolimits_{w\left(  \widetilde{Q}\right)
}B\\\text{(by Claim 5)}}}\det\underbrace{\left(  \operatorname*{sub}%
\nolimits_{w\left(  Q\right)  }^{w\left(  K\right)  }D\right)  }%
_{\substack{=I_{k}\\\text{(by Claim 6)}}}\ \ \ \ \ \ \ \ \ \ \left(
\text{since }\widetilde{\widetilde{Q}}=Q\right) \nonumber\\
&  =\left(  -1\right)  ^{\sum\widetilde{Q}+\sum\widetilde{K}}\det\left(
\operatorname*{rows}\nolimits_{w\left(  \widetilde{Q}\right)  }B\right)
\underbrace{\det\left(  I_{k}\right)  }_{=1}\nonumber\\
&  =\left(  -1\right)  ^{\sum\widetilde{Q}+\sum\widetilde{K}}\det\left(
\operatorname*{rows}\nolimits_{w\left(  \widetilde{Q}\right)  }B\right)  .
\label{sol.det.AB=0-proport.detD=}%
\end{align}

We know that $K$ and $\widetilde{K}$ are two subsets of $\left\{
1,2,\ldots,n\right\}  $ satisfying $\left\vert K\right\vert +\left\vert
\widetilde{K}\right\vert =n$. Also, $CD=\left(  g_{i,j}\right)  _{1\leq i\leq
n,\ 1\leq j\leq n}\in\mathbb{K}^{n\times n}$ is an $n\times n$-matrix such
that%
\[
\text{every }i\in K\text{ and }j\in\widetilde{K}\text{ satisfy }g_{i,j}=0
\]
(by Claim 7). Hence, Exercise \ref{exe.det.blocktria-twisted} \textbf{(b)}
(applied to $K$, $\widetilde{K}$, $CD$ and $g_{i,j}$ instead of $P$, $Q$, $A$
and $a_{i,j}$) yields%
\begin{align*}
\det\left(  CD\right)   &  =\left(  -1\right)  ^{\sum K+\sum
\widetilde{\widetilde{K}}}\det\left(  \operatorname*{sub}\nolimits_{w\left(
K\right)  }^{w\left(  \widetilde{\widetilde{K}}\right)  }\left(  CD\right)
\right)  \det\left(  \operatorname*{sub}\nolimits_{w\left(  \widetilde{K}%
\right)  }^{w\left(  \widetilde{K}\right)  }\left(  CD\right)  \right) \\
&  =\underbrace{\left(  -1\right)  ^{\sum K+\sum K}}%
_{\substack{=1\\\text{(since }\sum K+\sum K=2\sum K\\\text{is even)}}%
}\det\underbrace{\left(  \operatorname*{sub}\nolimits_{w\left(  K\right)
}^{w\left(  K\right)  }\left(  CD\right)  \right)  }%
_{\substack{=\operatorname*{cols}\nolimits_{w\left(  Q\right)  }A\\\text{(by
Claim 8)}}}\underbrace{\det\left(  \operatorname*{sub}\nolimits_{w\left(
\widetilde{K}\right)  }^{w\left(  \widetilde{K}\right)  }\left(  CD\right)
\right)  }_{\substack{=\operatorname*{rows}\nolimits_{w\left(  \widetilde{P}%
\right)  }B\\\text{(by Claim 9)}}}\\
&  \ \ \ \ \ \ \ \ \ \ \ \ \ \ \ \ \ \ \ \ \left(  \text{since }%
\widetilde{\widetilde{K}}=K\right) \\
&  =\det\left(  \operatorname*{cols}\nolimits_{w\left(  Q\right)  }A\right)
\det\left(  \operatorname*{rows}\nolimits_{w\left(  \widetilde{P}\right)
}B\right)  .
\end{align*}
Therefore,%
\begin{align*}
&  \det\left(  \operatorname*{cols}\nolimits_{w\left(  Q\right)  }A\right)
\det\left(  \operatorname*{rows}\nolimits_{w\left(  \widetilde{P}\right)
}B\right) \\
&  =\det\left(  CD\right)  =\det C\cdot\det D\ \ \ \ \ \ \ \ \ \ \left(
\text{by (\ref{sol.det.AB=0-proport.detCD=})}\right) \\
&  =\left(  -1\right)  ^{\sum\widetilde{K}+\sum\widetilde{P}}\det\left(
\operatorname*{cols}\nolimits_{w\left(  P\right)  }A\right)  \cdot\left(
-1\right)  ^{\sum\widetilde{Q}+\sum\widetilde{K}}\det\left(
\operatorname*{rows}\nolimits_{w\left(  \widetilde{Q}\right)  }B\right) \\
&  \ \ \ \ \ \ \ \ \ \ \ \ \ \ \ \ \ \ \ \ \left(  \text{here, we multiplied
the equalities (\ref{sol.det.AB=0-proport.detC=}) and
(\ref{sol.det.AB=0-proport.detD=})}\right) \\
&  =\left(  -1\right)  ^{\sum\widetilde{K}+\sum\widetilde{P}}\left(
-1\right)  ^{\sum\widetilde{Q}+\sum\widetilde{K}}\cdot\det\left(
\operatorname*{cols}\nolimits_{w\left(  P\right)  }A\right)  \det\left(
\operatorname*{rows}\nolimits_{w\left(  \widetilde{Q}\right)  }B\right)  .
\end{align*}
In view of%
\begin{align*}
&  \left(  -1\right)  ^{\sum\widetilde{K}+\sum\widetilde{P}}\left(  -1\right)
^{\sum\widetilde{Q}+\sum\widetilde{K}}\\
&  =\left(  -1\right)  ^{\left(  \sum\widetilde{K}+\sum\widetilde{P}\right)
+\left(  \sum\widetilde{Q}+\sum\widetilde{K}\right)  }\\
&  =\left(  -1\right)  ^{\sum\widetilde{P}+\sum\widetilde{Q}}%
\ \ \ \ \ \ \ \ \ \ \left(
\begin{array}
[c]{c}%
\text{since }\left(  \sum\widetilde{K}+\sum\widetilde{P}\right)  +\left(
\sum\widetilde{Q}+\sum\widetilde{K}\right) \\
=\sum\widetilde{P}+\sum\widetilde{Q}+2\sum\widetilde{K}\equiv\sum
\widetilde{P}+\sum\widetilde{Q}\operatorname{mod}2
\end{array}
\right)  ,
\end{align*}
we can rewrite this as%
\begin{align*}
&  \det\left(  \operatorname*{cols}\nolimits_{w\left(  Q\right)  }A\right)
\det\left(  \operatorname*{rows}\nolimits_{w\left(  \widetilde{P}\right)
}B\right) \\
&  =\left(  -1\right)  ^{\sum\widetilde{P}+\sum\widetilde{Q}}\det\left(
\operatorname*{cols}\nolimits_{w\left(  P\right)  }A\right)  \det\left(
\operatorname*{rows}\nolimits_{w\left(  \widetilde{Q}\right)  }B\right)  .
\end{align*}
Multiplying both sides of this equality by $\left(  -1\right)  ^{\sum
\widetilde{P}}$, we obtain%
\begin{align*}
&  \left(  -1\right)  ^{\sum\widetilde{P}}\det\left(  \operatorname*{cols}%
\nolimits_{w\left(  Q\right)  }A\right)  \cdot\det\left(  \operatorname*{rows}%
\nolimits_{w\left(  \widetilde{P}\right)  }B\right) \\
&  =\underbrace{\left(  -1\right)  ^{\sum\widetilde{P}}\left(  -1\right)
^{\sum\widetilde{P}+\sum\widetilde{Q}}}_{\substack{=\left(  -1\right)
^{\sum\widetilde{P}+\left(  \sum\widetilde{P}+\sum\widetilde{Q}\right)
}=\left(  -1\right)  ^{\sum\widetilde{Q}}\\\text{(since }\sum\widetilde{P}%
+\left(  \sum\widetilde{P}+\sum\widetilde{Q}\right)  =2\sum\widetilde{P}%
+\sum\widetilde{Q}\equiv\sum\widetilde{Q}\operatorname{mod}2\text{)}}%
}\det\left(  \operatorname*{cols}\nolimits_{w\left(  P\right)  }A\right)
\det\left(  \operatorname*{rows}\nolimits_{w\left(  \widetilde{Q}\right)
}B\right) \\
&  =\left(  -1\right)  ^{\sum\widetilde{Q}}\det\left(  \operatorname*{cols}%
\nolimits_{w\left(  P\right)  }A\right)  \cdot\det\left(  \operatorname*{rows}%
\nolimits_{w\left(  \widetilde{Q}\right)  }B\right)  .
\end{align*}
This solves Exercise \ref{exe.det.AB=0-proport}.
\end{proof}

\subsection{Solution to Exercise \ref{exe.cauchy-border}}

Our solution to Exercise \ref{exe.cauchy-border} will rely on the following
simple lemma:\footnote{Definition \ref{def.submatrix.minor} is being used
here.}

\begin{lemma}
\label{lem.sol.cauchy-border.lapxy}Let $n\in\mathbb{N}$. Let $A=\left(
a_{i,j}\right)  _{1\leq i\leq n,\ 1\leq j\leq n}$ be an $n\times n$-matrix.
Let $x_{1},x_{2},\ldots,x_{n}$ be $n$ elements of $\mathbb{K}$. Let
$y_{1},y_{2},\ldots,y_{n}$ be $n$ elements of $\mathbb{K}$. Then,%
\[
\sum_{i=1}^{n}\ \ \sum_{j=1}^{n}\left(  -1\right)  ^{i+j}a_{i,j}\left(
x_{i}+y_{j}\right)  \det\left(  A_{\sim i,\sim j}\right)  =\left(  \sum
_{k=1}^{n}x_{k}+\sum_{k=1}^{n}y_{k}\right)  \cdot\det A.
\]

\end{lemma}

We might refer to Lemma \ref{lem.sol.cauchy-border.lapxy} as \textquotedblleft
Laplace expansion along all rows and columns simultaneously\textquotedblright.
The following proof explains why:

\begin{proof}
[Proof of Lemma \ref{lem.sol.cauchy-border.lapxy}.]We have%
\begin{align*}
&  \sum_{i=1}^{n}\ \ \sum_{j=1}^{n}\left(  -1\right)  ^{i+j}a_{i,j}\left(
x_{i}+y_{j}\right)  \det\left(  A_{\sim i,\sim j}\right) \\
&  =\sum_{p=1}^{n}\ \ \sum_{q=1}^{n}\underbrace{\left(  -1\right)
^{p+q}a_{p,q}\left(  x_{p}+y_{q}\right)  \det\left(  A_{\sim p,\sim q}\right)
}_{=\left(  -1\right)  ^{p+q}a_{p,q}x_{p}\det\left(  A_{\sim p,\sim q}\right)
+\left(  -1\right)  ^{p+q}a_{p,q}y_{q}\det\left(  A_{\sim p,\sim q}\right)
}\\
&  \ \ \ \ \ \ \ \ \ \ \ \ \ \ \ \ \ \ \ \ \left(
\begin{array}
[c]{c}%
\text{here, we have renamed the}\\
\text{summation indices }i\text{ and }j\text{ as }p\text{ and }q
\end{array}
\right) \\
&  =\sum_{p=1}^{n}\ \ \sum_{q=1}^{n}\left(  \left(  -1\right)  ^{p+q}%
a_{p,q}x_{p}\det\left(  A_{\sim p,\sim q}\right)  +\left(  -1\right)
^{p+q}a_{p,q}y_{q}\det\left(  A_{\sim p,\sim q}\right)  \right) \\
&  =\sum_{p=1}^{n}\ \ \sum_{q=1}^{n}\left(  -1\right)  ^{p+q}a_{p,q}x_{p}%
\det\left(  A_{\sim p,\sim q}\right)  +\underbrace{\sum_{p=1}^{n}%
\ \ \sum_{q=1}^{n}}_{=\sum_{q=1}^{n}\ \ \sum_{p=1}^{n}}\left(  -1\right)
^{p+q}a_{p,q}y_{q}\det\left(  A_{\sim p,\sim q}\right) \\
&  =\sum_{p=1}^{n}\ \ \underbrace{\sum_{q=1}^{n}\left(  -1\right)
^{p+q}a_{p,q}x_{p}\det\left(  A_{\sim p,\sim q}\right)  }_{=x_{p}\sum
_{q=1}^{n}\left(  -1\right)  ^{p+q}a_{p,q}\det\left(  A_{\sim p,\sim
q}\right)  }+\sum_{q=1}^{n}\ \ \underbrace{\sum_{p=1}^{n}\left(  -1\right)
^{p+q}a_{p,q}y_{q}\det\left(  A_{\sim p,\sim q}\right)  }_{=y_{q}\sum
_{p=1}^{n}\left(  -1\right)  ^{p+q}a_{p,q}\det\left(  A_{\sim p,\sim
q}\right)  }\\
&  =\sum_{p=1}^{n}x_{p}\underbrace{\sum_{q=1}^{n}\left(  -1\right)
^{p+q}a_{p,q}\det\left(  A_{\sim p,\sim q}\right)  }_{\substack{=\det
A\\\text{(by Theorem \ref{thm.laplace.gen} \textbf{(a)})}}}+\sum_{q=1}%
^{n}y_{q}\underbrace{\sum_{p=1}^{n}\left(  -1\right)  ^{p+q}a_{p,q}\det\left(
A_{\sim p,\sim q}\right)  }_{\substack{=\det A\\\text{(by Theorem
\ref{thm.laplace.gen} \textbf{(b)})}}}\\
&  =\sum_{p=1}^{n}x_{p}\det A+\sum_{q=1}^{n}y_{q}\det A=\left(  \sum_{p=1}%
^{n}x_{p}\right)  \det A+\left(  \sum_{q=1}^{n}y_{q}\right)  \det A\\
&  =\left(  \sum_{p=1}^{n}x_{p}+\sum_{q=1}^{n}y_{q}\right)  \cdot\det
A=\left(  \sum_{k=1}^{n}x_{k}+\sum_{k=1}^{n}y_{k}\right)  \cdot\det A
\end{align*}
(here, we have renamed the summation indices $p$ and $q$ as $k$ and
$k$).\ This proves Lemma \ref{lem.sol.cauchy-border.lapxy}.
\end{proof}

We are now ready to solve Exercise \ref{exe.cauchy-border}:

\begin{proof}
[Solution to Exercise \ref{exe.cauchy-border}.]\textbf{(a)} The definition of
$C$ yields $C=\left(  \dfrac{1}{x_{i}+y_{j}}\right)  _{1\leq i\leq n,\ 1\leq
j\leq n}$. Hence, Lemma \ref{lem.sol.cauchy-border.lapxy} (applied to $A=C$
and $a_{i,j}=\dfrac{1}{x_{i}+y_{j}}$) yields%
\[
\sum_{i=1}^{n}\ \ \sum_{j=1}^{n}\left(  -1\right)  ^{i+j}\dfrac{1}{x_{i}%
+y_{j}}\left(  x_{i}+y_{j}\right)  \det\left(  C_{\sim i,\sim j}\right)
=\left(  \sum_{k=1}^{n}x_{k}+\sum_{k=1}^{n}y_{k}\right)  \cdot\det C.
\]
Comparing this with%
\[
\sum_{i=1}^{n}\ \ \sum_{j=1}^{n}\left(  -1\right)  ^{i+j}\underbrace{\dfrac
{1}{x_{i}+y_{j}}\left(  x_{i}+y_{j}\right)  }_{=1}\det\left(  C_{\sim i,\sim
j}\right)  =\sum_{i=1}^{n}\ \ \sum_{j=1}^{n}\left(  -1\right)  ^{i+j}%
\det\left(  C_{\sim i,\sim j}\right)  ,
\]
we obtain%
\begin{align}
&  \sum_{i=1}^{n}\ \ \sum_{j=1}^{n}\left(  -1\right)  ^{i+j}\det\left(
C_{\sim i,\sim j}\right) \nonumber\\
&  =\left(  \sum_{k=1}^{n}x_{k}+\sum_{k=1}^{n}y_{k}\right)  \cdot\det C.
\label{sol.cauchy-border.a.1}%
\end{align}

However, the definition of $\operatorname*{adj}C$ (see Definition
\ref{def.adj}) yields \newline$\operatorname*{adj}C=\left(  \left(  -1\right)
^{i+j}\det\left(  C_{\sim j,\sim i}\right)  \right)  _{1\leq i\leq n,\ 1\leq
j\leq n}$. Thus,
\begin{align}
&  \left(  \text{the sum of all entries of }\operatorname*{adj}C\right)
\nonumber\\
&  =\underbrace{\sum_{i=1}^{n}\ \ \sum_{j=1}^{n}}_{=\sum_{j=1}^{n}%
\ \ \sum_{i=1}^{n}}\underbrace{\left(  -1\right)  ^{i+j}}_{=\left(  -1\right)
^{j+i}}\det\left(  C_{\sim j,\sim i}\right) \label{sol.cauchy-border.a.2a}\\
&  =\sum_{j=1}^{n}\ \ \sum_{i=1}^{n}\left(  -1\right)  ^{j+i}\det\left(
C_{\sim j,\sim i}\right)  =\sum_{i=1}^{n}\ \ \sum_{j=1}^{n}\left(  -1\right)
^{i+j}\det\left(  C_{\sim i,\sim j}\right) \nonumber\\
&  \ \ \ \ \ \ \ \ \ \ \ \ \ \ \ \ \ \ \ \ \left(
\begin{array}
[c]{c}%
\text{here, we have renamed the}\\
\text{summation indices }i\text{ and }j\text{ as }j\text{ and }i
\end{array}
\right) \nonumber\\
&  =\left(  \sum_{k=1}^{n}x_{k}+\sum_{k=1}^{n}y_{k}\right)  \cdot\det
C\ \ \ \ \ \ \ \ \ \ \left(  \text{by (\ref{sol.cauchy-border.a.1})}\right)  .
\label{sol.cauchy-border.a.2}%
\end{align}
In view of%
\begin{align*}
\det C  &  =\det\left(  \left(  \dfrac{1}{x_{i}+y_{j}}\right)  _{1\leq i\leq
n,\ 1\leq j\leq n}\right)  \ \ \ \ \ \ \ \ \ \ \left(  \text{since }C=\left(
\dfrac{1}{x_{i}+y_{j}}\right)  _{1\leq i\leq n,\ 1\leq j\leq n}\right) \\
&  =\frac{\prod_{1\leq i<j\leq n}\left(  \left(  x_{i}-x_{j}\right)  \left(
y_{i}-y_{j}\right)  \right)  }{\prod_{\left(  i,j\right)  \in\left\{
1,2,\ldots,n\right\}  ^{2}}\left(  x_{i}+y_{j}\right)  }%
\ \ \ \ \ \ \ \ \ \ \left(  \text{by Exercise \ref{exe.cauchy-det}}\right)  ,
\end{align*}
we can rewrite this as%
\begin{align}
&  \left(  \text{the sum of all entries of }\operatorname*{adj}C\right)
\nonumber\\
&  =\left(  \sum_{k=1}^{n}x_{k}+\sum_{k=1}^{n}y_{k}\right)  \cdot\frac
{\prod_{1\leq i<j\leq n}\left(  \left(  x_{i}-x_{j}\right)  \left(
y_{i}-y_{j}\right)  \right)  }{\prod_{\left(  i,j\right)  \in\left\{
1,2,\ldots,n\right\}  ^{2}}\left(  x_{i}+y_{j}\right)  }.
\label{sol.cauchy-border.a.3}%
\end{align}
This solves Exercise \ref{exe.cauchy-border} \textbf{(a)}.

\textbf{(b)} Assume that the matrix $C$ is invertible. Then, Theorem
\ref{thm.matrices.inverses.square} \textbf{(a)} (applied to $A=C$) yields that
the matrix $C$ is invertible if and only if the element $\det C$ of
$\mathbb{K}$ is invertible (in $\mathbb{K}$). Hence, the element $\det C$ of
$\mathbb{K}$ is invertible (since the matrix $C$ is invertible). Therefore,
Theorem \ref{thm.matrices.inverses.square} \textbf{(b)} (applied to $A=C$)
yields that the inverse of $C$ is $C^{-1}=\dfrac{1}{\det C}\cdot
\operatorname*{adj}C$.

\begin{vershort}
In other words, each entry of the matrix $C^{-1}$ equals $\dfrac{1}{\det C}$
times the corresponding entry of $\operatorname*{adj}C$. Therefore,%
\begin{align*}
&  \left(  \text{the sum of all entries of }C^{-1}\right) \\
&  =\dfrac{1}{\det C}\cdot\left(  \text{the sum of all entries of
}\operatorname*{adj}C\right)  =\sum_{k=1}^{n}x_{k}+\sum_{k=1}^{n}%
y_{k}\ \ \ \ \ \ \ \ \ \ \left(  \text{by (\ref{sol.cauchy-border.a.2}%
)}\right)  .
\end{align*}

\end{vershort}

\begin{verlong}
Thus, for each $\left(  i,j\right)  \in\left\{  1,2,\ldots,n\right\}  ^{2}$,
we have%
\begin{align}
&  \left(  \text{the }\left(  i,j\right)  \text{-th entry of }C^{-1}\right)
\nonumber\\
&  =\left(  \text{the }\left(  i,j\right)  \text{-th entry of }\dfrac{1}{\det
C}\cdot\operatorname*{adj}C\right)  \ \ \ \ \ \ \ \ \ \ \left(  \text{since
}C^{-1}=\dfrac{1}{\det C}\cdot\operatorname*{adj}C\right) \nonumber\\
&  =\dfrac{1}{\det C}\cdot\left(  \text{the }\left(  i,j\right)  \text{-th
entry of }\operatorname*{adj}C\right)  . \label{sol.cauchy-border.b.ij}%
\end{align}

Now, $\operatorname*{adj}C$ is an $n\times n$-matrix. Hence,%
\begin{align}
&  \left(  \text{the sum of all entries of }\operatorname*{adj}C\right)
\nonumber\\
&  =\sum_{i=1}^{n}\ \ \sum_{j=1}^{n}\left(  \text{the }\left(  i,j\right)
\text{-th entry of }\operatorname*{adj}C\right)  .
\label{sol.cauchy-border.b.adjC}%
\end{align}

Also, $C^{-1}$ is an $n\times n$-matrix. Hence,%
\begin{align*}
&  \left(  \text{the sum of all entries of }C^{-1}\right) \\
&  =\sum_{i=1}^{n}\ \ \sum_{j=1}^{n}\underbrace{\left(  \text{the }\left(
i,j\right)  \text{-th entry of }C^{-1}\right)  }_{\substack{=\dfrac{1}{\det
C}\cdot\left(  \text{the }\left(  i,j\right)  \text{-th entry of
}\operatorname*{adj}C\right)  \\\text{(by (\ref{sol.cauchy-border.b.ij}))}}}\\
&  =\sum_{i=1}^{n}\ \ \sum_{j=1}^{n}\dfrac{1}{\det C}\cdot\left(  \text{the
}\left(  i,j\right)  \text{-th entry of }\operatorname*{adj}C\right) \\
&  =\dfrac{1}{\det C}\cdot\underbrace{\sum_{i=1}^{n}\ \ \sum_{j=1}^{n}\left(
\text{the }\left(  i,j\right)  \text{-th entry of }\operatorname*{adj}%
C\right)  }_{\substack{=\left(  \text{the sum of all entries of }%
\operatorname*{adj}C\right)  \\\text{(by (\ref{sol.cauchy-border.b.adjC}))}%
}}\\
&  =\dfrac{1}{\det C}\cdot\underbrace{\left(  \text{the sum of all entries of
}\operatorname*{adj}C\right)  }_{\substack{=\left(  \sum_{k=1}^{n}x_{k}%
+\sum_{k=1}^{n}y_{k}\right)  \cdot\det C\\\text{(by
(\ref{sol.cauchy-border.a.2}))}}}\\
&  =\dfrac{1}{\det C}\cdot\left(  \sum_{k=1}^{n}x_{k}+\sum_{k=1}^{n}%
y_{k}\right)  \cdot\det C=\sum_{k=1}^{n}x_{k}+\sum_{k=1}^{n}y_{k}.
\end{align*}

\end{verlong}

\noindent This solves Exercise \ref{exe.cauchy-border} \textbf{(b)}.

\textbf{(c)} We shall use the notation from Definition
\ref{def.sol.det.rk1upd.scalar}. The definition of $D$ yields $D=\left(
\begin{array}
[c]{cc}%
C & u\\
u^{T} & 0_{1\times1}%
\end{array}
\right)  $.

We have $0_{1\times1}=\left(
\begin{array}
[c]{c}%
0
\end{array}
\right)  $. Thus, Corollary \ref{cor.sol.det.bordered.3} (applied to $v=u^{T}$
and $h=0$ and $H=0_{1\times1}$ and $A=C$) yields
\[
\det\left(
\begin{array}
[c]{cc}%
C & u\\
u^{T} & 0_{1\times1}%
\end{array}
\right)  =\underbrace{0\det C}_{=0}-\operatorname*{ent}\left(  u^{T}\left(
\operatorname*{adj}C\right)  u\right)  =-\operatorname*{ent}\left(
u^{T}\left(  \operatorname*{adj}C\right)  u\right)  .
\]
Since $D=\left(
\begin{array}
[c]{cc}%
C & u\\
u^{T} & 0_{1\times1}%
\end{array}
\right)  $, we can rewrite this as
\begin{equation}
\det D=-\operatorname*{ent}\left(  u^{T}\left(  \operatorname*{adj}C\right)
u\right)  . \label{sol.cauchy-border.c.1}%
\end{equation}

\begin{vershort}
However, $u\in\mathbb{K}^{n}$ was defined to be the column vector whose all
$n$ entries equal $1$. Thus, $u=\left(  1,1,\ldots,1\right)  ^{T}$ and
$u^{T}=\left(  1,1,\ldots,1\right)  $.
\end{vershort}

\begin{verlong}
However, $u\in\mathbb{K}^{n}$ was defined to be the column vector whose all
$n$ entries equal $1$. Thus, $u=\left(
\begin{array}
[c]{c}%
1\\
1\\
\vdots\\
1
\end{array}
\right)  =\left(  1,1,\ldots,1\right)  ^{T}$, and therefore
\begin{align*}
u^{T}  &  =\left(
\begin{array}
[c]{c}%
1\\
1\\
\vdots\\
1
\end{array}
\right)  ^{T}\ \ \ \ \ \ \ \ \ \ \left(  \text{since }u=\left(
\begin{array}
[c]{c}%
1\\
1\\
\vdots\\
1
\end{array}
\right)  \right) \\
&  =\left(  1,1,\ldots,1\right)  .
\end{align*}

\end{verlong}

Therefore, Corollary \ref{cor.sol.det.rk1upd.claim.ent(v(adjA)u)} (applied to
$u_{i}=1$, $v_{i}=1$, $v=u^{T}$ and $A=C$) yields%
\begin{align*}
\operatorname*{ent}\left(  u^{T}\left(  \operatorname*{adj}C\right)  u\right)
&  =\sum_{i=1}^{n}\ \ \sum_{j=1}^{n}\left(  -1\right)  ^{i+j}%
\underbrace{1\cdot1}_{=1}\det\left(  C_{\sim j,\sim i}\right) \\
&  =\sum_{i=1}^{n}\ \ \sum_{j=1}^{n}\left(  -1\right)  ^{i+j}\det\left(
C_{\sim j,\sim i}\right) \\
&  =\left(  \text{the sum of all entries of }\operatorname*{adj}C\right)
\ \ \ \ \ \ \ \ \ \ \left(  \text{by (\ref{sol.cauchy-border.a.2a})}\right) \\
&  =\left(  \sum_{k=1}^{n}x_{k}+\sum_{k=1}^{n}y_{k}\right)  \cdot\frac
{\prod_{1\leq i<j\leq n}\left(  \left(  x_{i}-x_{j}\right)  \left(
y_{i}-y_{j}\right)  \right)  }{\prod_{\left(  i,j\right)  \in\left\{
1,2,\ldots,n\right\}  ^{2}}\left(  x_{i}+y_{j}\right)  }%
\ \ \ \ \ \ \ \ \ \ \left(  \text{by (\ref{sol.cauchy-border.a.3})}\right)  .
\end{align*}
In light of this, we can rewrite (\ref{sol.cauchy-border.c.1}) as%
\[
\det D=-\left(  \sum_{k=1}^{n}x_{k}+\sum_{k=1}^{n}y_{k}\right)  \cdot
\frac{\prod_{1\leq i<j\leq n}\left(  \left(  x_{i}-x_{j}\right)  \left(
y_{i}-y_{j}\right)  \right)  }{\prod_{\left(  i,j\right)  \in\left\{
1,2,\ldots,n\right\}  ^{2}}\left(  x_{i}+y_{j}\right)  }.
\]
This solves Exercise \ref{exe.cauchy-border} \textbf{(c)}.
\end{proof}

\section{Appendix: Old citations}

Some places (articles, internet postings, etc.) cite older versions of these
notes. Most of the time, the numbering of theorems, propositions, formulas and
even chapters has shifted since these older versions were written; thus, these
citations no longer point to the correct location of the results they want to
cite. The following tables list the results cited in all citations known to
me, and their locations in the current version of these notes.%

\[
\text{%
\begin{tabular}
[c]{|l|l|}\hline\hline
item cited & its current location\\\hline\hline
& \\\hline
Additional exercise 23 in the version of 19 May 2016 & Exercise
\ref{exe.det.rk1upd}\\\hline
Corollary 5.44 in the version of 29 May 2016 & Corollary
\ref{cor.laplace.pre.col}\\\hline
Corollary 5.85 in the version of 29 May 2016 & Corollary \ref{cor.adj.kernel}%
\\\hline
Exercise 41 in the version of 4 September 2016 & Exercise \ref{exe.det.2diags}%
\\\hline
Lemma 5.49 in the version of 2016-09-04 & Lemma \ref{lem.vander-det.lem-rearr}%
\\\hline
Exercise 36 in the version of 6 October 2016 & Exercise
\ref{exe.vander-det.xi+yj}\\\hline
Proposition 2.25 in the version of 2016-12-22 & Proposition
\ref{prop.vandermonde.consequences}\\\hline
Corollary 7.10 in the version of 2016-12-22 & Corollary
\ref{cor.AoPS262752.X=n}\\\hline
Theorem 5.157 in the version of 15 February 2017 & Theorem \ref{thm.det(A+B)}%
\\\hline
Corollary 5.161 in the version of 25 May 2017 & Corollary \ref{cor.det(A+X)}%
\\\hline
Additional exercise 22 in the version of 25 May 2017 & Exercise
\ref{addexe.jacobi-complement}\\\hline
Theorem 5.32 in the version of 25 May 2017 & Theorem \ref{thm.cauchy-binet}%
\\\hline
Theorem 5.157 in the version of 25 May 2017 & Theorem \ref{thm.det(A+B)}%
\\\hline
Exercise 44 in the version of 28 October 2017 & Exercise \ref{exe.det.2diags}%
\\\hline
Example 5.28 in the version of 7 November 2017 & Example \ref{exam.xi+yj}%
\\\hline
Exercise 5.16 in the version of 18 November 2017 & Exercise
\ref{exe.vander-det.xi+yj}\\\hline
\end{tabular}
}%
\]%
\[
\text{%
\begin{tabular}
[c]{|l|l|}\hline\hline
item cited & its current location\\\hline\hline
& \\\hline
Section 5.5 in the version of 14 March 2018 & Section \ref{sect.det.CB}%
\\\hline
Exercise 5.33 in the version of 21 March 2018 & Exercise \ref{exe.adj(AB)}%
\\\hline
Corollary 5.163 in the version of 21 March 2018 & Corollary \ref{cor.det(A+X)}%
\\\hline
Theorem 2.27 in the version of 26 April 2018 & Theorem
\ref{thm.vandermonde.XY}\\\hline
Corollary 5.161 in the version of 26 April 2018 & Corollary \ref{cor.det(A+D)}%
\\\hline
Corollary 7.63 in the version of 26 April 2018 & Corollary
\ref{cor.ps2.2.S.-1alt}\\\hline
Corollary 7.150 in the version of 26 April 2018 & Corollary
\ref{cor.adj(AB).cauchy-binet-general}\\\hline
Exercise 5.21 in the version of 14 September 2018 & Exercise
\ref{exe.det.a1a2anx}\\\hline
Remark 2.27 in the version of 4 October 2018 & Remark \ref{rmk.ind.abstract}%
\\\hline
\S 3.1 in the version of 4 October 2018 & Section \ref{sect.binom.def}\\\hline
Theorem 3.29 in the version of 4 October 2018 & Theorem
\ref{thm.vandermonde.rat}\\\hline
Theorem 3.44 in the version of 4 October 2018 & Theorem \ref{thm.pie.jordan}%
\\\hline
Corollary 8.63 in the version of 4 October 2018 & Corollary
\ref{cor.ps2.2.S.-1}\\\hline
Chapter 6 in the version of 18 October 2018 & Chapter \ref{chp.det}\\\hline
\end{tabular}
}%
\]

\end{document}